\documentclass[12pt]{amsart}
\usepackage[a4paper]{anysize}
\usepackage{amsthm}
\usepackage{amssymb}
\usepackage{times}
\usepackage{mathrsfs}
\usepackage{mathtools}
\usepackage{bbm}
\usepackage[unicode]{hyperref}

\include{xy}
\xyoption{v2}
\xyoption{all}
\xyoption{2cell}
\UseAllTwocells

    \newtheorem{theorem}{Theorem}[subsection]
    \newtheorem{proposition}[theorem]{Proposition}
    \newtheorem{lemma}[theorem]{Lemma}
    \newtheorem{corollary}[theorem]{Corollary}
    
    \theoremstyle{definition}
    \newtheorem{example}[theorem]{Example}
    \newtheorem{definition}[theorem]{Definition}
    \newtheorem{remark}[theorem]{Remark}
    \theoremstyle{remark}
    \newtheorem{claim}[theorem]{Claim}

    \def\set{\setcounter{equation}
             {\value{theorem}}\addtocounter{theorem}{1}}
    \numberwithin{equation}{subsection}
    \def\sset{\setcounter{subsubsection}
             {\value{theorem}}\addtocounter{theorem}{1}}

\makeatletter
\def\@tocline#1#2#3#4#5#6#7{\relax
   \ifnum #1>\c@tocdepth 
   \else
     \par \addpenalty\@secpenalty\addvspace{#2}%
     \begingroup \hyphenpenalty\@M
     \@ifempty{#4}{%
       \@tempdima\csname r@tocindent\number#1\endcsname\relax
     }{%
       \@tempdima#4\relax
     }%
     \parindent\z@ \leftskip#3\relax \advance\leftskip\@tempdima\relax
     \rightskip\@pnumwidth plus4em \parfillskip-\@pnumwidth
     #5\leavevmode\hskip-\@tempdima #6\nobreak\relax
     \ifnum#1<0\hfill\else\dotfill\fi\hbox to\@pnumwidth{\@tocpagenum{#7}}\par
     \nobreak
     \endgroup
   \fi}
\makeatother
    
\def\A{{\mathbb A}}
\def\bA{{\mathbf A}}
\def\cA{{{\mathscr A}}}

\def\B{{\mathbb B}}
\def\bB{{\mathbf B}}
\def\cB{{{\mathscr B}}}
\def\C{{\mathbb C}}

\def\cC{{{\mathscr C}}}
\def\sA{\mathsf{A}}
\def\sB{\mathsf{B}}
\def\sC{\mathsf{C}}
\def\cD{{{\mathscr D}}}

\def\fD{{\mathfrak D}}
\def\sD{\mathsf{D}}
\def\E{{\mathbb E}}
\def\cE{{{\mathscr E}}}
\def\bE{{\mathbf E}}
\def\sE{\mathsf{E}}
\def\F{{\mathbb F}}
\def\cF{{{\mathscr F}}}

\def\sF{\mathsf{F}}
\def\sG{\mathsf{G}}
\def\G{{\mathbb G}}
\def\cG{{{\mathscr G}}}

\def\cH{{{\mathscr H}}}

\def\bI{{\mathbb I}}
\def\cI{{{\mathscr I}}}
\def\cJ{{{\mathscr J}}}
\def\K{{\mathbb K}}
\def\cK{{{\mathscr K}}}
\def\cL{{{\mathscr L}}}
\def\L{{\mathbb L}}
\def\bL{{\mathbf L}}
\def\sL{\mathsf{L}}
\def\cM{{{\mathscr M}}}
\def\fM{{\mathfrak M}}
\def\sM{\mathsf{M}}
\def\N{{\mathbb N}}
\def\sN{\mathsf{N}}
\def\cN{{{\mathscr N}}}
\def\cO{{{\mathscr O}}}

\def\cP{{{\mathscr P}}}
\def\fp{{\mathfrak p}}
\def\P{{\mathbb P}}
\def\bP{{\mathbf P}}
\def\fP{{\mathfrak P}}
\def\sP{{\mathsf P}}
\def\sQ{{\mathsf Q}}
\def\bQ{{\mathbf Q}}
\def\cQ{{{\mathscr Q}}}
\def\fq{{\mathfrak q}}
\def\Q{{\mathbb Q}}
\def\R{{\mathbb R}}
\def\cR{{{\mathscr R}}}
\def\cS{{{\mathscr S}}}
\def\SS{{\mathbb S}}
\def\cT{{{\mathscr T}}}
\def\T{{\mathbb T}}
\def\bT{{\mathbf T}}
\def\sT{{{\mathsf T}}}

\def\tt{\mathrm{t}}
\def\sS{\mathsf{S}}
\def\sT{\mathsf{T}}
\def\sU{\mathsf{U}}
\def\sV{\mathsf{V}}
\def\cX{{{\mathscr X}}}
\def\sX{{{\mathsf X}}}
\def\sY{{{\mathsf Y}}}
\def\cY{{{\mathscr Y}}}
\def\Z{{\mathbb Z}}
\def\bZ{{\mathbf Z}}
\def\cZ{{{\mathscr Z}}}
\def\sZ{{{\mathsf Z}}}
\def\cW{{{\mathscr W}}}
\def\W{{\mathbb W}}
\def\V{{\mathbb V}}
\def\cU{{{\mathscr U}}}
\def\cV{{{\mathscr V}}}
\def\cK{{{\mathscr K}}}

\makeatletter
\newcommand{\superimpose}[2]{%
  {\ooalign{$#1\@firstoftwo#2$\cr\hfil$#1\@secondoftwo#2$\hfil\cr}}}
\makeatother

\def\boxvert{\:\mathpalette\superimpose{{\square}{\rule[0.1mm]{.3pt}{2.7mm}}}\:}

\newcommand{\colim}
{\mathop\mathrm{colim}}

\newcommand{\Lim}
{\mathop\mathrm{Lim}}

\newcommand{\Colim}
{\mathop\mathrm{Colim}}

\newcommand{\Pslim}[1]
{\underset{#1}{\mathrm{2\text{-}lim}\,}}

\newcommand{\Pscolim}[1]
{\underset{#1}{\mathrm{2\text{-}colim}\,}}

\newcommand{\derotimes}
{\overset{\mathbf{L}}{\otimes}}

\newcommand{\otimesint}
{\overset{{}_\mathrm{int}}{\otimes}}

\newcommand{\ellotimes}
{\overset{\ell}{\otimes}}

\renewcommand{\labelenumi}{(\roman{enumi})}
\renewcommand{\labelenumii}{(\alph{enumii})}

\newenvironment{pfclaim}[1][$\Diamond$]
{\def\claimQED{{#1}}\noindent {\em Proof of the claim. }}
{\leavevmode\unskip\penalty9999 \hbox{}\nobreak\hfill
    \quad\hbox{\claimQED}{\smallskip}}

\newenvironment{acknowledgement}{\vskip .5cm
\setlength{\baselineskip}{0mm}
\noindent \footnotesize\rmfamily
{\em Acknowledgements\ \ \ }}

\def\Can{\mathrm{Can}}
\def\Tr{\mathrm{Tr}}
\def\Str{\mathrm{Str}}
\def\tr{\mathrm{tr}}
\def\trunc{\mathrm{trunc}}
\def\gp{\mathrm{gp}}
\def\gr{\mathrm{gr}}

\def\alt{\mathrm{alt}}
\def\red{\mathrm{red}}
\def\sep{\mathrm{sep}}
\def\Ann{\mathrm{Ann}}
\def\Pol{\mathrm{Pol}}
\def\perf{\mathrm{perf}}

\def\depth{\mathrm{depth}}

\def\Cone{\mathrm{Cone}}
\def\Cont{\mathrm{Cont}}
\def\Gal{\mathrm{Gal}}
\def\Spec{\mathrm{Spec}}
\def\Spa{\mathrm{Spa}}
\def\Spv{\mathrm{Spv}}
\def\Proj{\mathrm{Proj}}
\def\Bl{\mathrm{Bl}}
\def\satBl{\mathrm{sat.Bl}}
\def\Alhom{\mathrm{alHom}}
\def\Hom{\mathrm{Hom}}
\def\Fun{\mathrm{Fun}}
\def\bFun{\mathsf{Fun}}

\def\sCart{\mathsf{Cart}}
\def\cCart{\mathscr{C}\!art}
\def\Fib{\mathsf{Fib}}
\def\LFib{\mathsf{LFib}}
\def\fib{\mathsf{fib}}
\def\cFib{\mathscr{F}ib}
\def\cspFib{sp\mathscr{F}ib}
\def\sFib{\mathsf{spFib}}
\def\lex{\mathsf{lex}}
\def\sLink{\mathsf{Link}}
\def\wLink{\mathsf{wLink}}
\def\Split{\mathsf{Split}}
\def\Pt{\mathrm{Pt}}
\def\PreStack{\mathsf{PreStack}}

\def\Stack{\mathsf{Stack}}
\def\LStack{\mathsf{LStack}}
\def\StGpd{\mathsf{StGpd}}

\def\St{\mathsf{St}}
\def\sMorph{\mathsf{Morph}}
\def\rMorph{\mathrm{Morph}}

\def\sPsFun{\mathsf{PsFun}}
\def\uniPsFun{\mathsf{uniPsFun}}
\def\stPsFun{\mathsf{stPsFun}}
\def\sPsNat{\mathsf{PsNat}}
\def\Hot{\mathsf{Hot}}
\def\rHot{\mathrm{Hot}}
\def\Der{\mathrm{Der}}

\def\Ext{\mathrm{Ext}}
\def\hExt{{\mathbb{E}\mbox{\sf xt}}}

\def\Exal{\mathrm{Exal}}
\def\Exaltop{\mathrm{Exaltop}}
\def\sExal{\mathsf{Exal}}
\def\snilExal{\mathsf{nilExal}}
\def\sExaltop{\mathsf{Exaltop}}
\def\snilExaltop{\mathsf{nilExaltop}}
\def\Tors{\mathsf{Tors}}

\def\Tor{\mathrm{Tor}}
\def\End{\mathrm{End}}
\def\cEnd{\mathscr{E}\!nd}

\def\Ind{\mathrm{Ind}}
\def\Res{\mathrm{Res}}
\def\Spf{\mathrm{Spf}}
\def\Coker{\mathrm{Coker}}
\def\Ker{\mathrm{Ker}}
\def\Img{\mathrm{Im}}
\def\Aut{\mathrm{Aut}}
\def\Alg{{{\text{-}}\mathsf{Alg}}}
\def\dga{{{\text{-}}\mathsf{dgAlg}}}
\def\dgalt{{{\text{-}}\mathsf{dgAlt}}}
\def\mAlg{{{\text{-}}\mathbf{m.Alg}}}
\def\Sym{\mathrm{Sym}}

\def\AlgMod{\mathsf{Alg.Mod}}
\def\dAlgMod{\text{-}\mathsf{dAlg.Mod}}

\def\hor{\mathrm{hor}}
\def\horf{\mathrm{hor.fl}}
\def\fl{\mathrm{fl}}
\def\Desc{\mathsf{Desc}}
\def\Et{\text{-}\mathbf{\acute{E}t}}
\def\Zar{\mathrm{Zar}}
\def\et{\mathrm{\acute{e}t}}

\def\wEt{\text{-}\mathbf{w.\acute{E}t}}

\def\bCat{\mathbf{Cat}}
\def\bGraph{\mathbf{Graph}}
\def\bNbd{\mathbf{Nbd}}
\def\bLift{\mathbf{Lift}}
\def\Mod{\text{-}\mathbf{Mod}}
\def\CM{\text{-}\mathsf{CM}}
\def\IntMnd{\mathbf{Int.Mnd}}
\def\SatMnd{\mathbf{Sat.Mnd}}
\def\Fan{\mathbf{Fan}}
\def\Mnd{\mathbf{Mnd}}

\def\Set{\mathbf{Set}}
\def\Site{\mathbf{Site}}
\def\wSite{\mathbf{wSite}}
\def\totSite{\mathsf{totSite}}
\def\Topos{\mathbf{Topos}}
\def\chara{\mathrm{char}}
\def\intg{\mathrm{int}}
\def\nil{\mathrm{nil}}

\def\rk{\mathrm{rk}}
\def\one{\mathbf{1}}
\def\bone{\mathbbm{1}}
\def\zero{\mathbf{0}}
\def\bar#1{\overline{#1}}
\def\hat{\widehat} \def\tilde{\widetilde}
\def\fm{\mathfrak m}
\def\fr{\mathfrak r}
\def\eps{\varepsilon}
\def\sGamma{\mathsf\Gamma}
\def\sSpec{\mathsf{Spec}}

\def\fSet{\mathbf{f.Set}}
\renewcommand\emptyset{\varnothing}
\renewcommand\phi{\varphi}
\renewcommand\theta{\vartheta}
\def\rad{\mathrm{rad}}
\def\Max{\mathrm{Max}}
\def\Min{\mathrm{Min}}
\def\sPhi{\mathsf{\Phi}}
\def\sPsi{\mathsf{\Psi}}
\def\Grp{\mathbf{Grp}}
\def\Gpd{\mathsf{Gpd}}
\def\Poset{\mathbf{POSet}}

\def\cDer{\mathscr{D}\!er}
\def\cHom{\mathscr{H}\!om}

\def\cTor{\mathscr{T}\!or}
\def\ba{\mathbf{a}}
\def\bj{\mathbf{j}}
\def\be{\mathbf{e}}
\def\bee{\boldsymbol{e}}
\def\bff{\mathbf{f}}
\def\bg{\mathbf{g}}

\def\bk{\mathbf{k}}
\def\bek{\boldsymbol{k}}
\def\bep{\boldsymbol{p}}
\def\bx{\mathbf{x}}
\def\bK{\mathbf{K}}
\def\br{\mathbf{r}}
\def\bt{\mathbf{t}}

\def\fn{\mathfrak{n}}
\def\Supp{\mathrm{Supp}}
\def\length{\mathrm{length}}
\def\sh{\mathrm{sh}}
\def\rmst{\mathrm{st}}
\def\coh{\mathrm{coh}}
\def\cohs{{\mathrm{coh},\{s\}}}
\def\qcoh{\mathrm{qcoh}}
\def\sqcoh{\mathsf{qcoh}}
\def\qcC{\mathsf{q}\cC}
\def\qcqs{\mathsf{qcqs}}
\def\he{\mathrm{h}}
\def\loc{\mathrm{loc}}
\def\isom{\stackrel{\sim}{\to}}
\def\hgt{\mathrm{ht}}
\def\hdim{\mathrm{proj.dim}}

\def\Ass{\mathrm{Ass}}
\def\Ob{\mathrm{Ob}}
\def\Div{\mathrm{Div}}
\def\cDiv{\mathscr{D}iv}
\def\cohDiv{\mathrm{coh.Div}}
\def\bDiv{\mathbf{Div}}
\def\Equal{\mathrm{Equal}}
\def\Equiv{\mathsf{Equiv}}
\def\Frac{\mathrm{Frac}}
\def\bRflx{\text{-}\mathbf{Rflx}}

\def\bprelog{\mathbf{pre\text{-}log}}
\def\prelog{\mathrm{pre\text{-}log}}
\def\blog{\mathbf{log}}

\def\logMod{\mathbf{Mod.pre\text{-}log}}
\def\bcohlog{\mathbf{coh.log}}
\def\bintlog{\mathbf{int.log}}
\def\bsatlog{\mathbf{sat.log}}
\def\bqcohlog{\mathbf{qcoh.log}}
\def\bqflog{\mathbf{qf.log}}
\def\bqfslog{\mathbf{qfs.log}}
\def\bfslog{\mathbf{fs.log}}
\def\bflog{\mathbf{f.log}}
\def\bPic{\mathbf{Pic}}
\def\Pic{\mathrm{Pic}}
\def\grPic{\mathbf{gr.Pic}}
\def\sdet{\mathsf{det}}
\def\Dec{\mathrm{D\acute{e}c}\,}
\def\Fil{\mathrm{Fil}}
\def\sFil{\mathsf{Fil}}
\def\SpSeq{\mathsf{Sp.Seq}}
\def\Lef{\mathrm{Lef}}
\def\Tate{\mathsf{Tate}}
\def\Top{\mathbf{Top}}
\def\sTop{\mathsf{Top}}
\def\TopMon{\mathsf{TopMon}}
\def\TopAlg{\mathsf{TopAlg}}
\def\pregrTopAlg{\mathsf{pre}\tdu\mathsf{gr.TopAlg}}
\def\grTopAlg{\mathsf{gr.TopAlg}}
\def\TopMod{\mathsf{TopMod}}
\def\TopHens{\mathsf{TopHens}}
\def\rmTop{\mathrm{Top}}
\def\topo{\mathrm{top}}
\def\Tot{\mathrm{Tot}}
\def\reg{\mathrm{reg}}
\def\sm{\mathrm{sm}}
\def\lfft{\mathrm{lfft}}
\def\sat{\mathrm{sat}}
\def\Syz{\mathrm{Syz}}
\def\fQ{{\mathfrak Q}}
\def\fU{{\mathfrak U}}

\def\fX{{\mathfrak X}}
\def\fZ{{\mathfrak Z}}
\def\sa{\mathsf{a}}
\def\sbb{\mathsf{b}}
\def\sc{\mathsf{c}}
\def\sd{\mathsf{d}}
\def\se{\mathsf{e}}
\def\sff{\mathsf{f}}
\def\slf{\mathsf{lf}}

\def\shh{\mathsf{h}}
\def\si{\mathsf{i}}
\def\sj{\mathsf{j}}
\def\ssm{\mathsf{m}}
\def\sp{\mathsf{p}}
\def\sq{\mathsf{q}}
\def\st{\mathsf{t}}
\def\su{\mathsf{u}}
\def\ss{\mathsf{s}}
\def\sx{\mathsf{x}}
\def\sy{\mathsf{y}}
\def\sR{\mathsf{R}}
\def\bsb{{\{s\}}}

\def\Alt{\mathsf{Alt}}
\def\Galois{\mathbf{Galois}}
\def\Perf{\mathsf{Perf}}
\def\Sch{\mathsf{Sch}}
\def\bCov{\mathsf{Cov}}
\def\bTame{\mathbf{Tame}}
\def\sMA{\mathsf{MA}}

\def\bdelta{{\boldsymbol{\delta}}}

\def\bDelta{\mathbf{\Delta}}
\def\bLambda{\mathbf{\Lambda}}
\def\sDelta{\mathsf{\Delta}}

\def\bd{\boldsymbol{d}}
\def\bm{{\mathbf{m}}}

\def\bp{{\boldsymbol{p}}}
\def\bmu{\boldsymbol{\mu}}
\def\etab{\boldsymbol{\eta}}
\def\beps{\boldsymbol{\eps}}
\def\bphi{\boldsymbol{\phi}}
\def\bpsi{\boldsymbol{\psi}}
\def\blambda{{\boldsymbol{\lambda}}}
\def\bomega{\boldsymbol{\omega}}
\def\bsigma{\boldsymbol{\sigma}}
\def\bOmega{\mathbf{\Omega}}

\def\alphaenu{\renewcommand{\labelenumi}{(\alph{enumi})}}
\def\romanenu{\renewcommand{\labelenumi}{(\roman{enumi})}}
\def\romanenuii{\renewcommand{\labelenumii}{(\roman{enumii})}}
\def\addenu{\addtocounter{enumi}{1}}
\def\addenuii{\addtocounter{enumii}{1}}

\def\La{\langle}
\def\Ra{\rangle}
\def\tdu{\text{-}}
\def\mapsfrom{\leftarrow\!\shortmid}
\def\squig{\rightsquigarrow}
\def\cosk{\mathrm{cosk}}
\def\sk{\mathrm{sk}}
\def\Sk{\mathrm{Sk}}
\def\Sh{\mathsf{Sh}}
\def\AW{\mathsf{AW}}
\def\ev{\mathrm{ev}}
\def\sev{\mathsf{ev}}
\def\ZR{\mathrm{ZR}}
\def\sign{\mathrm{sign}}

\makeindex

\begin{document}
\title{Foundations for almost ring theory}

\author{Ofer Gabber}
\author{Lorenzo Ramero}

\maketitle

\centerline\today

\vskip 1cm

\centerline{\textsf{\slshape Release 7.5}}

\vskip 2cm

    Ofer Gabber

    I.H.E.S.

    Le Bois-Marie

    35, route de Chartres

    F-91440 Bures-sur-Yvette

    {\em e-mail address:} \href{mailto:gabber@ihes.fr}{gabber@ihes.fr}

\vskip .5cm

    Lorenzo Ramero

    Universit{\'e} de Lille I

    Laboratoire de Math{\'e}matiques

    F-59655 Villeneuve d'\!Ascq C\'edex

{\em e-mail address:} \href{mailto:ramero@math.univ-lille1.fr}
{ramero@math.univ-lille1.fr}

{\em web page:} \url{http://math.univ-lille1.fr/~ramero}

\vskip 2cm

\begin{acknowledgement} We thank Niels Borne, Kestutis Cesnavicius,
Michel Emsalem, Pierre-Yves Gaillard, Lutz Geissler, W.-P. Heidorn,
Fabrice Orgogozo, Olaf Schn\"urer and Peter Scholze for pointing out
some mistakes in earlier drafts, and for useful suggestions.
\end{acknowledgement}

\newpage

\tableofcontents


\hfill{\small\rm It is not incumbent upon you to complete the work,}

\hfill{\small\rm but neither are you at liberty to desist from it.}

\hfill{\small\rm (Avot 2:21)}

\setcounter{section}{-1}
\section{Introduction}

Both the focus of this monograph and its subject matter
have evolved considerably in the last few years.
On the one hand, the insistence on making the text
self-contained (aside from a reduced canon of basic
references, which should ideally contain only EGA and
some parts of Bourbaki's {\em \'El\'ements}), has resulted
in a rather weighty mass of material of independent
interest, that is applied to, but is completely
separate from almost ring theory, and whose relationship
to $p$-adic Hodge theory is thus even more indirect.
Rather than stemming from a well-thought-out plan, this
part is the outcome of a haphazard process, lumbering
between alternating phases of accretion and consolidation,
with new topics piled up as dictated by need, or occasionally
by whim, when we just branched out from the main flow to
pursue a certain line of thought to its logical destination.
Nevertheless, a few themes have spontaneously emerged,
around which the originally amorphous magma has been
able to settle, to the point where by now a distinct
shape is finally discernible, and it is perhaps time to
pause and take stock of its broad outlines.

Now then, we may distinguish :

$\bullet$\ \
First of all, a rather thorough exposition of the foundations
of logarithmic algebraic geometry, comprising chapters
\ref{chap_monoids}, \ref{chap_logschemes} and \ref{chap_etale-cov}.
Inevitably, our treatment owes a lot to the works of Kato
and his school : our contribution is foremost that of
gathering and tidying up the subject, which until now
was scattered in a disparate number of research articles,
many of them still unpublished and even unfinished.
A closer scrutiny would also reveal a few technical
innovations that we hope will become standard issue
of the working algebraic log-geometer : we may mention
the systematic use of pointed monoids and pointed modules,
the projective fan associated with a graded monoid, or a
definition of $\alpha$-flatness for log structures which
refines and generalizes an older notion of ``toric flatness''.
Furthermore, we took the occasion to repair a few small
(and not so small) mistakes and inaccuracies that we
detected in the literature.

$\bullet$\ \
Two other chapters are dedicated to local cohomology and
Grothendieck's duality theory. Early on, the emphasis here
was on generalizations : especially, we were interested
in removing from the theory the pervasive noetherian
assumptions, to pave the way for our recasting of Faltings's
almost purity in the framework of valuation theory.
Applications of local cohomology to non-archimedean
analytic geometry furnished another influential motivation,
though one that has remained, so far, hidden from view.
More recently, the noetherian aspects have also become
relevant to our project, and this latest release contains
a detailed account of the most important properties of
noetherian rings endowed with a dualizing complex. The
latter, in turn, could be dealt with satisfactorily only
after a thorough revisitation of the general theory of the
dualizing complex, so that our chapters \ref{chap_local-coh}
and \ref{chap_duality} can also be regarded as complementary
to Conrad's book \cite{Con} (dedicated to the trace morphism
and the deeper aspects of duality) : totaling our respective
efforts, it should eventually become possible to bypass
entirely Hartshorne's notes \cite{Har} which, as is well
known, are wanting in many ways.

$\bullet$\ \
Chapters \ref{chap_hom-algebra} and \ref{chap_comm-algebra}
present (for the time being, anyway) a looser structure :
a miscellanea of self-contained units devoted to more or
less independent topics. However, there is at least one
thread running through several sections, and whose stretch
can be traced all the way back to the earliest beginnings of
almost ring theory; it connects sections \ref{sec_simplicial}
and \ref{sec_simplicial-sets} -- on simplicial homotopy theory --
to a section \ref{sec_homotopy} dedicated to homotopical
algebra, then on to sections \ref{sec_rings-regular} and
\ref{sec_excellent-rings}, which make extensive use of the
cotangent complex to derive important characterizations
of regular and excellent rings, including an up-to-date
presentation of classical results due to Andr\'e, extracted
from his monograph \cite{An}, and from his paper \cite{AnII}
on localization of formal smoothness. This homotopical
algebraic thread resurfaces again in section
\ref{sec_non-flat}, but there we are already squarely into
almost ring theory proper.
\vskip.5cm
On the other hand, two recent notable developments are
compelling a revision of our understanding of almost ring
theory itself, and of its situation within commutative
algebra and algebraic geometry at large :

$\bullet$\ \
The first is Scholze's PhD thesis \cite{Scho} on
{\em perfectoid spaces}, that contains both a maximal
generalization of the almost purity theorem, and a major
simplification of its proof, based on his ``tilting''
technique (and completely different from Faltings's).
However, the range of Scholze's theory transcends the
domain of $p$-adic Hodge theory (which was not even his
original motivation) : to drive the point home, his thesis
concludes with a clever application to the long standing
weight monodromy conjecture, thus affording the unusual
spectacle of a tool which was fashioned out of purely
$p$-adic concerns, and ends up playing a crucial role
in the solution of a purely $\ell$-adic problem.

$\bullet$\ \
The second spectacular development is Yves Andr\'e's proof
of the direct summand conjecture (\cite{YAn}); the latter is
a deceptively simple assertion that has been a central problem
in commutative algebra for the last thirty years : it asserts
that every finite injective ring homomorphism $f:A\to B$
from a regular local ring $A$, admits a $A$-linear splitting.
The relevance of almost purity to this question was first
surmised by Paul Roberts in 2001 (after a talk by the second
author at the University of Utah), and has been widely
advertised by him ever since. Andr\'e's solution uses
perfectoid techniques, and builds on earlier work by Bhargav
Bhatt, who in \cite{Bha} proved the conjecture in the case
where $A$ is essentially smooth over a mixed characteristic
discrete valuation ring and $f\otimes_\Z\Q$ is \'etale outside
a relative normal crossings divisor of $\Spec\,A$. Moreover,
Bhatt has subsequently simplified some of Andr\'e's arguments
and shown how the same method yields a more general
``derived version'' of the conjecture, for proper schemes over
any regular ring : see \cite{Bha2}.

We see then, that almost ring theory has emancipated itself
from its former ancillary role in the exclusive service of
$p$-adic Hodge theory, and is now elbowing out a niche in
the wider ecosystem of algebraic geometry.

The present release completes the project announced
in the introduction of the 6th release : 

$\bullet$\ \ 
First we introduce a class of topological rings that generalize
the perfectoid rings of \cite{Scho}; it is very easy to say what
a (generalized) perfectoid $\F_p$-algebra is : namely, it is
just a perfect and complete topological $\F_p$-algebra whose
topology is linear, defined by an ideal of finite type.
The general definition is somewhat more involved, but we
prove the following characterization. For any perfectoid
$\F_p$-algebra $E$, we consider the ring of Witt vectors
$W(E)$, and we endow it with a natural topology, induced
from that of $E$; then every perfectoid ring is a
topological quotient of the form $A:=W(E)/\underline aW(E)$,
for such a suitable $E$, and where $\underline a\in W(E)$
is what we call a {\em distinguished element} : see
definition \ref{def_distinguished}. Moreover, just as in
Scholze's work, the perfectoid $\F_p$-algebra $E$ can
be recovered from $A$ via a {\em tilting functor} that
establishes, more precisely, an equivalence between the
category of all perfectoid rings and that of pairs $(E,\cI)$
consisting of a perfectoid $\F_p$-algebra $E$ and a principal
ideal $\cI\subset W(E)$ generated by a distinguished element
(as it is well known, this construction is rooted in Fontaine
and Winterberger's theory of the {\em field of norms}).
The distinguished ideal $\cI$ represents an extra parameter
that remains hidden in Scholze's original approach : the reason
is that he fixes from the start a base perfectoid field $K$,
thereby implicitly fixing as well a distinguished element
$\underline a$ in the ring of Witt vectors of the tilt of
$K$, and then every perfectoid ring in his work is supposed
to be a $K$-algebra, which -- from our viewpoint -- amounts
to restricting to perfectoid rings whose associated distinguished
ideal is generated by $\underline a$. Having thus removed the
parameter $\cI$, he can then also do away with Witt vectors
altogether, and the inverse to the tilting construction is
obtained in \cite{Scho} via a more abstract deformation theoretic
argument. This route is precluded to us, so we rely instead
on direct and rather concrete Witt vectors calculations.
A similar strategy has been proposed in \cite{Ked-Liu}, and
our viewpoint can indeed be described fairly as an interpolation
of those of Scholze and Kedlaya-Liu, though we only studied
\cite{Scho} in detail.

$\bullet$\ \
The first three sections of chapter \ref{chap_perfectoid}
are devoted to exploring this new class of perfectoid rings
and its manifold remarkable features. The rest of the chapter
then merges the theory of perfectoid rings with Huber's adic
spaces, to forge the perfectoid spaces that are the main tool
for our proof of almost purity, whose most general form is
given by theorem \ref{th-alm-purity-form-perfectoid} and
applies to {\em formal perfectoid rings}, {\em i.e.} to
topological rings whose completion is perfectoid. The proof
proceeds via several preliminary reductions : first, to the
case of a {\em perfectoid quasi-affinoid ring}, covered by
theorem \ref{th_a-purity-for-perfectoids}, then -- by exploiting
the local geometry of adic spaces -- to the case of a perfectoid
valuation ring, which was treated already in our monograph
\cite{Ga-Ra}. What enables here this localization argument
is a basic feature of the \'etale topology of arbitrary adic
spaces : the fibred category of finite \'etale coverings of
the affinoid subsets of an adic space is a stack. The latter
result is in turn a special case of our theorem
\ref{th_only-for-ind-finite}.

$\bullet$\ \
We also include a detailed treatment of the foundations of
the theory of adic spaces, that essentially follows \cite{Hu1},
but contains some modest improvement : notably, the
systematic use of {\em analytically noetherian rings} (borrowed
from \cite{Fu-Ga-Ka}) allows us to unify the two classes of
topological rings that Huber dealt with separately in
his work (the strongly noetherian rings and the f-adic
rings with a noetherian ring of definition). We also
point out a henselian variant of the structure sheaf
that is available on the adic spectrum of any f-adic
ring, with no restriction whatsoever.

$\bullet$\ \
The last chapter proposes a few applications :
in section \ref{sec_model-alg} we introduce a class of
{\em model algebras} over any rank one valuation ring $K^+$
of mixed characteristic $(0,p)$, and we show that when $K^+$
is deeply ramified, such algebras are formal perfectoid
rings for their $p$-adic topology; hence the theory of
chapter \ref{chap_perfectoid} immediately yields an almost
purity theorem for model algebras. Likewise, section
\ref{sec_regular-pure} proves an almost purity theorem for
certain very ramified towers of log-regular rings; again,
after some preliminary reductions, the proof amounts to the
observation that the inductive limits of such towers are formal
perfectoid for their $p$-adic topology. These instances of almost
purity were already contained in a previous draft of our work
(Release 6), where they were proven by an extension of Faltings's
method, that relied on deep results from logarithmic algebraic
geometry, and also entailed the construction of certain
{\em normalized lengths} for torsion modules over model
algebras, and respectively over the rings occurring in
section \ref{sec_regular-pure}. Neither of these two ingredients
intervenes any longer in the new proofs; however, we have found
worthwhile to explain how model algebras arise from suitable
very ramified towers of log-smooth $K^+$-algebras, and we have
also retained the construction of normalized lengths for model
algebras and for limits of towers log-regular rings, since they
are sufficiently interesting in their own right, and might be
useful for other applications (normalized lengths for
torsion $K^+$-modules are exploited in \cite{Scho2}).

Section \ref{sec_dir-summand-conj} contains our account of
Andr\'e's work on the direct summand conjecture, which we
generalize to the case of a finite injective ring homomorphism
$A\to B$ where $A$ is {\em log-regular} : see theorem
\ref{th_log-direct-summand-conj}, whose proof relies
on a refinement of Andr\'e's {\em perfectoid Abhyankar's lemma},
which is the main result of section \ref{sec_perf-Abhyankar}.
In a future release we shall present a corresponding logarithmic
generalization of Bhatt's derived version of the direct summand
conjecture.

Section \ref{sec_diagonal-idempotents} considers a finite
injective map $f:A\to B$ of noetherian rings, where $A$ is
again log-regular; if $f$ is \'etale, $B\otimes_AB$ has a
canonical {\em diagonal idempotent} $e_f$, so called because
its support is the open and closed diagonal in
$\Spec\,B\times_{\Spec\,A}\Spec\,B$. When $f$ is only generically
\'etale, $e_f$ is only well defined on some localization
$B\otimes_AB[a^{-1}]$ (for some $a\in A$ whose zero locus is
nowhere dense in $\Spec\,A$), and the problem we address,
is to exhibit an explicit $d\in B\otimes_AB$ such that $de_f$
is {\em integral}, {\em i.e.} lies in the image of $B\otimes_AB$.
This question has an easy answer when $B$ is a projective
$A$-module : indeed, in this case we have a well defined trace
map $B\to A$, whence a {\em different ideal} $\cD\subset B$,
and a standard calculation shows that $(b\otimes 1)\cdot e_f$
is integral for every $b\in\cD$. If $B$ is not projective, the
question seems to be much more difficult : our solution makes
an essential use of perfectoid techniques, and especially the
perfectoid Abhyankar's lemma of section
\ref{sec_perf-Abhyankar} : see theorem
\ref{th_diagonal-idempotent} and remark \ref{rem_Hulot-demissionne}.

The last section contains for now only some preliminary
construction that shall be used, in a future release, to
prove the existence and weak functoriality of {\em big
Cohen-Macaulay algebras}, extending previous work by Andr\'e
and others.



\section{Basic category theory}
The purpose of this chapter is to fix some notation
that shall stand throughout this work, and to collect,
for ease of reference, a few well known generalities on
categories and functors that are frequently used.
Our main reference on general nonsense is the treatise
\cite{Bor}, and another good reference is the more recent
\cite{Ka-Sch}.

Sooner or later, any honest discussion of categories and
topoi gets tangled up with some foundational issues revolving
around the manipulation of large sets. For this reason, to be
able to move on solid ground, it is essential to select
from the outset a definite set-theoretical framework (among
the several currently available), and stick to it unwaveringly.

Thus, {\em throughout this work we will accept the so-called
Zermelo-Fraenkel system of axioms for set theory}. (In this
version of set theory, everything is a set, and there is
no primitive notion of class, in contrast to other
axiomatisations.)

Additionally, following \cite[Exp.I, \S0]{SGA4-1}, we shall
assume that, for every set $S$, there exists a {\em universe\/}
$\sV$ such that $S\in\sV$. (For the notion of universe, the
reader may also see \cite[\S1.1]{Bor}.)

Throughout this chapter, we fix some universe $\sU$ such
that $\N\in\sU$ (where $\N$ is the set of natural numbers;
the latter condition is required, in order to be able to
perform some standard set-theoretical operations without
leaving $\sU$). A set $S$ is {\em $\sU$-small} (resp.
{\em essentially $\sU$-small}), if $S\in\sU$ (resp. if $S$
has the cardinality of a $\sU$-small set). If the context is
not ambiguous, we shall just write small, instead of $\sU$-small.

\subsection{Categories, functors and natural transformations}
\label{sec_label-for-morph}
A {\em category\/} $\cC$ is the datum of a set $\Ob(\cC)$ of
{\em objects\/} and, for every $A,B\in\Ob(\cC)$, a set of
{\em morphisms} from $A$ to $B$, denoted :
$$
\Hom_\cC(A,B)
$$
and as usual we write $f:A\to B$ to signify $f\in\Hom_\cC(A,B)$.
Furthermore, we set
$$
\rMorph(\cC):=\{(A,B,f)~|~A,B\in\Ob(\cC),\ f\in\Hom_\cC(A,B)\}.
$$
For any $\underline f:=(A,B,f)\in\rMorph(\cC)$, the object $A$ is
called the {\em source} of $\underline f$, and $B$ is the
{\em target} of $\underline f$. We also often use the notation
$$
\End_\cC(A):=\Hom_\cC(A,A)
$$
and the elements of $\End_\cC(A)$ are called the {\em endomorphisms}
of $A$ in $\cC$. We say that a pair of elements
$(\underline f,\underline g)$ of $\rMorph(\cA)$ is {\em composable}
if the target of $\underline f$ equals the source of $\underline g$.
Moreover, for every $A,B,C\in\Ob(\cC)$ we have a {\em composition law}
$$
\Hom_\cC(A,B)\times\Hom_\cC(B,C)\to\Hom_\cC(A,C)
\quad :\quad
(f,g)\mapsto g\circ f
$$
fulfilling the following two standard axioms :
\begin{itemize}
\item
For every $A\in\Ob(\cC)$ there exists an {\em identity endomorphism}
$\one_A$ of $A$, such that
$$
\one_A\circ f=f
\quad
g\circ\one_A=g
\qquad
\text{for every $B,C\in\Ob(\cC)$ and every $f:B\to A$ and $g:A\to C$}.
$$
\item
The composition law is {\em associative}, {\em i.e.} we have
$$
(h\circ g)\circ f=h\circ(g\circ f)
$$
for every $A,B,C,D\in\Ob(\cC)$ and every $f:A\to B$, $g:B\to C$
and $h:C\to D$.
\end{itemize}
Clearly, it follows that $(\End_\cC(A),\circ,\one_A)$ is
a monoid, and we get a group :
$$
\Aut_\cC(A)\subset\End_\cC(A)
$$
of invertible endomorphisms, {\em i.e.} the {\em automorphisms}
of the object $A$.

\sset\subsubsection{}\label{subsec_opposing}
We say that the category $\cC$ is {\em $\sU$-small} (or just
{\em small}), if both $\Ob(\cC)$ and $\rMorph(\cC)$ are small
sets. We say that $\cC$ {\em has small $\Hom$-sets} if
$\Hom_\cC(A,B)\in\sU$ for every $A,B\in\Ob(\cC)$.

A {\em subcategory} of $\cC$ is a category $\cB$ with
$\Ob(\cB)\subset\Ob(\cC)$ and $\rMorph(\cB)\subset\rMorph(\cC)$.

The {\em opposite category\/} $\cC^o$ is the category with
$\Ob(\cC^o)=\Ob(\cC)$, and such that :
$$
\Hom_{\cC^o}(A,B):=\Hom_\cC(B,A) \qquad\text{for every
$A,B\in\Ob(\cC)$}
$$
(with composition law induced by that of $\cC$, in the obvious way).
Given $A\in\Ob(\cC)$, sometimes we denote by $A^o$ the same object,
viewed as an element of $\Ob(\cC^o)$; likewise, given a morphism
$f:A\to B$ in $\cC$, we write $f^o$ for the corresponding morphism
$B^o\to A^o$ in $\cC^o$.

\sset\subsubsection{}\label{subsec_well-power}
A morphism $f:A\to B$ in $\cC$ is said to be a {\em monomorphism\/}
if the induced map
$$
\Hom_\cC(X,f):\Hom_\cC(X,A)\to\Hom_\cC(X,B)
\qquad
g\mapsto f\circ g
$$
is injective, for every $X\in\Ob(\cC)$. Dually, we say that $f$ is
an {\em epimorphism\/} if $f^o$ is a monomorphism in $\cC^o$.
Also, $f$ is an {\em isomorphism} if there exists a morphism
$g:B\to A$ such that $g\circ f=\one_A$ and $f\circ g=\one_B$.
Obviously, an isomorphism is both a monomorphism and an epimorphism.
The converse does not necessarily hold, in an arbitrary category.

Two monomorphisms $f:A\to B$ and $f':A'\to B$ are
{\em equivalent\/}, if there exists an isomorphism $h:A\to A'$ such
that $f=f'\circ h$. A {\em subobject\/} of $B$ is defined as an
equivalence class of monomorphisms $A\to B$. Dually, a
{\em quotient\/} of $B$ is a subobject of $B^o$ in $\cC^o$.

One says that $\cC$ is {\em well-powered} if, for every $A\in\Ob(\cC)$,
the set :
$$
\mathrm{Sub}(A)
$$
of all subobjects of $A$ is essentially small. Dually,
$\cC$ is {\em co-well-powered}, if $\cC^o$ is well-powered.

\sset\subsubsection{}
Let $\cA$ and $\cB$ be any two categories; a {\em functor}
$F:\cA\to\cB$ is a pair of maps
$$
\Ob(\cA)\to\Ob(\cB)
\qquad
\rMorph(\cA)\to\rMorph(\cB)
$$
both denoted also by $F$, such that
\begin{itemize}
\item
$F$ assigns to any morphism $f:A\to A'$ in $\cA$, a morphism
$Ff:FA\to FA'$ in $\cB$
\item
$F\one_A=\one_{FA}$ for every $A\in\Ob(\cA)$
\item
$F(g\circ f)=Fg\circ Ff$ for every $A,A',A''\in\Ob(\cA)$
and every pair of morphisms $f:A\to A'$, $g:A'\to A''$ in $\cA$.
\end{itemize}
If $F:\cA\to\cB$ and $G:\cB\to\cC$ are any two functors,
we get a composition
$$
G\circ F:\cA\to\cC
$$
which is the functor whose maps on objects and morphisms
are the compositions of the respective maps for $F$ and $G$.
We denote by
$$
\Fun(\cA,\cB)
$$
the set of all functors $\cA\to\cB$. Moreover, any such $F$
induces a functor $F^o:\cA^o\to\cB^o$ with $F^oA^o:=(FA)^o$
and $F^of^o:=(Ff)^o$ for every $A\in\Ob(\cA)$ and every
$f\in\rMorph(\cA)$.

\begin{definition}\label{def_equivalence}
Let $F:\cA\to\cB$ be a functor.

(i)\ \
We say that $F$ is {\em faithful} (resp. {\em full}, resp.
{\em fully faithful}), if for every $A,A'\in\Ob(\cA)$ it
induces injective (resp. surjective, resp. bijective) maps :
$$
\Hom_\cA(A,A')\to\Hom_\cB(FA,FA') \quad :\quad f\mapsto Ff.
$$

(ii)\ \
We say that $F$ {\em reflects monomorphisms} (resp. {\em reflects
epimorphisms}, resp. {\em is conservative}) if the following holds.
For every morphism $f:A\to A'$ in $\cA$, if the morphism $Ff$ of
$\cB$ is a monomorphism (resp. epimorphism, resp. isomorphism),
then the same holds for $f$.

(iii)\ \
If $\cA$ is a subcategory of $\cB$, and $F$ is the natural
inclusion functor, then $F$ is obviously faithful, and we
say that $\cA$ is a {\em full subcategory} of $\cB$, if $F$
is fully faithful.

(iv)\ \
The {\em essential image of $F$} is the full subcategory
of $\cB$ whose objects are the objects of $\cB$ that are
isomorphic to an object of the form $FA$, for some
$A\in\Ob(\cA)$.
We say that $F$ is {\em essentially surjective\/} if its
essential image is $\cB$.

(v)\ \
We say that $F$ is an {\em equivalence}, if it is fully
faithful and essentially surjective.
\end{definition}

\begin{remark}\label{rem_i-faithful}
For later use, it is convenient to introduce the notion
of {\em $n$-faithful functor}, for all integers $n\leq 2$.
Namely : if $n<0$, every functor is $n$-faithful; a functor
$F:\cA\to\cB$ (between any two categories $\cA$ and $\cB$)
is $0$-faithful, if it is faithful; $F$ is $1$-faithful,
if it is fully faithful; finally, we say that $F$ is
$2$-faithful, if it is an equivalence.
\end{remark}

\begin{example}\label{ex_universe}
(i)\ \
The collection of all small categories, together with the
functors between them, forms a category
$$
\sU\tdu\bCat.
$$
Unless we have to deal with more than one universe, we shall
usually omit the prefix $\sU$, and write just $\bCat$. It is
easily seen that $\bCat$ is a category with small $\Hom$-sets.

(ii)\ \
The category of all small sets shall be denoted $\sU\tdu\Set$
or just $\Set$, if there is no need to emphasize the chosen
universe. There is a natural fully faithful embedding :
$$
\Set\to\bCat.
$$
Indeed, to any set $S$ one may assign its {\em discrete category\/}
also denoted $S$, {\em i.e.} the unique category such that $\Ob(S)=S$
and $\rMorph(S)=\{(s,s,\one_s)~|~s\in S\}$. If $S$ and $S'$ are
two discrete categories, the datum of a functor $S\to S'$ is clearly
the same as a map of sets $\Ob(S)\to\Ob(S')$.

Notice also the natural functor
$$
\Ob:\bCat\to\Set
\qquad
\cC\mapsto\Ob(\cC)
$$
that assigns to each functor $F:\cC\to\cD$ the underlying
map $\Ob(\cC)\to\Ob(\cD)$ : $C\mapsto FC$.

(iii)\ \
Recall that a {\em preordered set} is a pair $(I,\leq)$
consisting of a set $I$ and a binary relation $\leq$ on
$I$ which is reflexive and transitive. In this case, we
also say that $\leq$ is a {\em preordering on $I$}.
We say that $(I,\leq)$ is a {\em partially ordered set},
if $\leq$ is also {\em antisymmetric}, {\em i.e.} if
we have
$$
(x\leq y\ \text{and}\ y\leq x)\Rightarrow x=y
\qquad
\text{for every $x,y\in I$}.
$$
We say that $(I,\leq)$ is a {\em totally ordered set}, if
it is partially ordered and any two elements are comparable,
{\em i.e.} for every $x,y\in I$ we have either $x\leq y$ or
$y\leq x$.  An {\em order-preserving map} $f:(I,\leq)\to(J,\leq)$
between preordered sets is a mapping $f:I\to J$ such that
$$
x\leq y\ \Rightarrow\ f(x)\leq f(y)
\qquad
\text{for every $x,y\in I$}.
$$
We denote by $\mathbf{Preorder}$ (resp. $\Poset$) the
category of small preordered (resp. partially ordered)
sets, with morphisms given by the order-preserving maps.
To any preordered set $(I,\leq)$ one assigns a category
whose set of objects is $I$, and whose morphisms are
given as follows. For every $i,j\in I$, the set of
morphisms $i\to j$ contains exactly one element when
$i\leq j$, and is empty otherwise. Clearly, this rule
defines a fully faithful functor
$$
\mathbf{Preorder}\to\bCat.
$$
Notice that if a category $C$ lies in the essential image
of this functor, then the same holds for $C^o$. Indeed, if
$C$ corresponds to the preordered set $(I,\leq)$, then
$C^o$ corresponds to the preordered set $(I^o,\leq)$
with $I^o:=I$ and $x\leq y$ in $I^o$ if and only if
$y\leq x$ in $I$, for every $x,y\in I$. Clearly $(I,\leq)$
is a partially ordered set if and only if the same holds
for $(I^o,\leq)$.
\end{example}

\sset\subsubsection{}
Let $\cA$, $\cB$ be two categories, $F,G:\cA\to\cB$ two
functors. A {\em natural transformation
\set\begin{equation}\label{eq_Rightarrow}
\alpha:F\Rightarrow G
\end{equation}
from $F$ to $G$\/} is a family of morphisms
$(\alpha_A:FA\to GA~|~A\in\Ob(\cA))$ of $\cB$ such that,
for every morphism $f:A\to B$ in $\cA$, the diagram :
\set\begin{equation}\label{eq_natural-transform}
{\diagram FA \ar[r]^-{\alpha_A} \ar[d]_{Ff} & GA \ar[d]^{Gf} \\
           FB \ar[r]^-{\alpha_B} & GB
\enddiagram}
\end{equation}
commutes. If $\alpha_A$ is an isomorphism for every $A\in\Ob(\cA)$,
we say that $\alpha$ is a {\em natural isomorphism\/} of functors.
For instance, the rule that assigns to any object $A$
the identity morphism $\one_{FA}:FA\to FA$, defines a natural
isomorphism $\one_F:F\Rightarrow F$. A natural transformation
\eqref{eq_Rightarrow} is also indicated by a diagram of the
type :
$$
\xymatrix{
\cA \rtwocell^F_G{\alpha} & \cB.
}$$

\sset\subsubsection{}\label{subsec_Godem-prod}
The natural transformations between functors $\cA\to\cB$
can be composed; namely, if $\alpha:F\Rightarrow G$ and
$\beta:G\Rightarrow H$ are two such transformations, we
obtain a natural transformation
$$
\beta\odot\alpha:F\Rightarrow H
\qquad\text{by the rule :}\qquad
A\mapsto\beta_A\circ\alpha_A
\qquad
\text{for every $A\in\Ob(\cA)$}.
$$
With this composition, $\Fun(\cA,\cB)$ is the set of objects of
a category which we shall denote 
$$
\bFun(\cA,\cB).
$$
There is also a second composition law for natural transformations :
if $\cC$ is another category, and we have a diagram of functors
and natural transformations
$$
\xymatrix{ \cA \rtwocell^F_G{\alpha} &
\cB \rtwocell^{F'}_{G'}{\alpha'} & \cC 
}$$
we get a natural transformation
$$
\alpha'*\alpha:F'\circ F\Rightarrow G'\circ G
\quad :\quad
A\mapsto \alpha'_{GA}\circ F'(\alpha_A)=G'(\alpha_A)\circ\alpha'_{FA}
\qquad
\text{for every $A\in\Ob(\cA)$}
$$
called the {\em Godement product\/} of $\alpha$ and $\alpha'$
(\cite[Prop.1.3.4]{Bor}).
Especially, if $H:\cB\to\cC$ (resp. $H:\cC\to\cA$) is any
functor, we write $H*\alpha$ (resp. $\alpha*H$) instead of
$\one_H*\alpha$ (resp. $\alpha*\one_H$).

Both composition laws are associative, {\em i.e.},
if we have additional natural transformations
$$
\xymatrix{
\cA\rtwocell^H_K{\gamma} & \cB &
\cC\rtwocell^{F''}_{G''}{\ \ \alpha''} & \cD
}$$
then we get the identities
$$
\gamma\odot(\beta\odot\alpha)=(\gamma\odot\beta)\odot\alpha
\qquad
\alpha''*(\alpha'*\alpha)=(\alpha''*\alpha')*\alpha.
$$
Moreover, the composition laws are related as follows.
Suppose that $\cA$, $\cB$ and $\cC$ are three categories,
and we have a diagram of six functors and four natural
transformations :
$$
\diagram
\cA\rruppertwocell<9>^{F_1}{\ \ \alpha_1}
\rrlowertwocell<-9>_{H_1}{\ \ \beta_1}
\ar[rr]|{G_1} & &
\cB\rruppertwocell<9>^{F_2}{\ \ \alpha_2}
\rrlowertwocell<-9>_{H_2}{\ \ \beta_2}
\ar[rr]|{G_2} & & \cC.
\enddiagram
$$
Then we have the identity :
\set\begin{equation}\label{eq_Godement-functor}
(\beta_2*\beta_1)\odot(\alpha_2*\alpha_1)=
(\beta_2\odot\alpha_2)*(\beta_1\odot\alpha_1).
\end{equation}
The proofs are left as exercises for the reader (see
\cite[Prop.1.3.5]{Bor}).

\begin{remark}
(i)\ \
In the situation of \eqref{subsec_Godem-prod}, if $\cA$ and
$\cB$ are small categories, the same holds for $\bFun(\cA,\cB)$.

(ii)\ \
Also, if $\cA$ is small, $\cB$ has small $\Hom$-sets and
$\Ob(\cB)\subset\sU$, then $\bFun(\cA,\cB)$ has small
$\Hom$-sets, and $\Ob(\bFun(\cA,\cB))\subset\sU$.

(iii)\ \
Assertion (ii) depends on our choices on how to encode
arbitrary maps of sets : according to our (implicit)
convention, {\em a map of sets $f:S\to S'$ is the graph
$\Gamma(f)\subset S\times S'$}. This does not agree,
{\em e.g.} with the definition found in Bourbaki's
treatise \cite{Bou-Ens}, where such a map $f$ is the
triple $(S,S',\Gamma(f))$. With Bourbaki's convention,
assertion (ii) fails. Other references are not so explicit
about their choices for encoding maps, but for instance the
SGA4 treatise (\cite{SGA4-1}, \cite{SGA4-2}, \cite{SGA4-3})
appears to follow Bourbaki's conventions, in view of
\cite[Exp.I, Rem.1.1.2]{SGA4-1}, which states that
$\Fun(\cA,\cB)$ is not necessarily a subset of $\sU$, and
$\bFun(\cA,\cB)$ does not necessarily have small $\Hom$-sets,
under the assumptions of (ii). On the other hand, under the
same assumptions, it is stated in \cite[Ch.II, Prop.1]{Ga} that
$\bFun(\cA,\cB)$ has small $\Hom$-sets, so the set-theoretic
conventions of the latter are not compatible with those of SGA4.
\end{remark}

\sset\subsubsection{Adjoint pairs of functors}
\label{subsec_adj-pair}
Let $\cA$ and $\cB$ be two categories, $F:\cA\to\cB$
and $G:\cB\to\cA$ two functors. We say that $G$ is
{\em left adjoint\/} to $F$ if there exist bijections
$$
\vartheta_{\!A,B}:\Hom_\cA(GB,A)\isom\Hom_\cB(B,FA)
\qquad
\text{for every $A\in\Ob(\cA)$ and $B\in\Ob(\cB)$}
$$
and these bijections are natural in both $A$ and $B$, {\em i.e.}
$$
\theta_{A'B'}(g\circ f\circ Gh)\!=\!Fg\circ\theta_{AB}(f)\circ h
\quad
\text{for all morphisms $GB\!\xrightarrow{f}\!A\!\xrightarrow{g}\!A'$
in $\cA$ and $B'\!\xrightarrow{h}\!B$ in $\cB$}.
$$
Then one also says that $F$ is {\em right adjoint\/} to $G$,
that $(G,F)$ is an {\em adjoint pair of functors}, and
that $\vartheta_{\bullet\bullet}$ is an {\em adjunction}
for the pair $(G,F)$.

Especially, to any object $B$ of $\cB$ (resp. $A$ of $\cA$),
the adjunction $\vartheta_{\bullet\bullet}$ assigns a morphism
$\theta_{GB,B}(\one_{GB}):B\to FGB$ (resp.
$\theta_{\!A,FA}^{-1}(\one_{FA}):GFA\to A$), whence a natural
transformation
\set\begin{equation}\label{eq_triangular}
\eta:\one_\cB\Rightarrow F\circ G
\qquad
\text{(resp. $\eps:G\circ F\Rightarrow\one_\cA$)}
\end{equation}
called the {\em unit\/} (resp. {\em counit}) of the adjunction.
The naturality of $\eta$ follows from the calculation:
$$
FG(f)\circ\eta_B=FG(f)\circ\theta_{GB,B}(\one_{GB})=
\theta_{GB,B'}(Gf\circ\one_{GB})=\theta_{GB,B'}(\one_{GB'}\circ Gf)
=\eta_{B'}\circ f
$$
for every morphism $f:B\to B'$ in $\cB$. A similar computation
shows the naturality of $\eps$. The naturality of
$\theta_{\bullet\bullet}$ implies that we have commutative diagrams
$$
\xymatrix{
B \ar[r]^-{\eta_B} \ar[rd]_{\theta_{A,B}(f)} &
FGB \ar[d]^{Ff} & 
GB \ar[r]^-{Gg} \ar[rd]_{\theta^{-1}_{A,B}(g)} & GFA
\ar[d]^{\eps_A} \\
& FA & & A
}$$
for every morphism $f:GB\to A$ in $\cA$ and $g:B\to FA$
in $\cB$. Taking $f=\eps_A$ and $g=\eta_B$, we see that
the unit and counit are related by the so-called
{\em triangular identities} expressed by the commutative
diagrams :
$$
\xymatrix{
F \ar@{=>}[r]^-{\eta*F} \ar@{=}[dr]_{\one_F} &
FGF \ar@{=>}[d]^{F*\eps} &
G \ar@{=>}[r]^-{G*\eta} \ar@{=}[dr]_{\one_G} &
GFG \ar@{=>}[d]^{\eps*G} \\
& F & & G.
}$$
Conversely, we have (\cite[Th.3.1.5]{Bor} or
\cite[Prop.1.5.4]{Ka-Sch}) :

\begin{proposition}\label{prop_triangular-identities}
Let $F:\cA\to\cB$ and $G:\cB\to\cA$ be two functors.
\begin{enumerate}
\item
If $\eta$, $\eps$ are natural transformations as in
\eqref{eq_triangular}, fulfilling the triangular
identities of \eqref{subsec_adj-pair}, then there
is a unique adjunction for the pair $(G,\!F)$, with
unit $\eta$ and counit $\eps$.
\item
Suppose that $(G,F)$ is an adjoint pair, and $\eta$ is the
unit (resp. $\eps$ is the counit) of an adjunction for
$(G,F)$, then there exists a unique natural transformation
$\eps$ (resp. $\eta$) as in \eqref{eq_triangular}, fulfilling
the triangular identities of \eqref{subsec_adj-pair}.
\end{enumerate}
\end{proposition}
\begin{proof}(i): Let $A\in\Ob(\cA)$ and $B\in\Ob(\cB)$ be
any two objects; in light of the discussion of
\eqref{subsec_adj-pair}, we see that the sought natural
bijection $\vartheta_{A,B}$ must be given by the rule :
$$
f\mapsto Ff\circ\eta_B
\qquad
\text{for every $f\in\Hom_\cA(GB,A)$}
$$
and its inverse must be the mapping
$$
g\mapsto\eps_A\circ Gg
\qquad
\text{for every $g\in\Hom_\cB(B,FA)$}.
$$
So we come down to checking that the triangular identities
imply that these rules do induce mutually inverse bijections
on the respective $\Hom$-sets. We leave the verification
to the reader.

(ii) is clear from the explicit construction of
$\vartheta_{\bullet\bullet}$ in (i).
\end{proof}

\begin{example}\label{ex_poset-quotient}
To every preordered set $(\cF,\leq)$ we may attach its
{\em partially ordered quotient}
$$
(\cF\!/\!\sim,\leq)
$$
where $\sim$ denotes the equivalence relation such that
$x\sim y$ if and only if $x\leq y$ and $y\leq x$, for
every $x,y\in\cF$. The ordering on $\cF\!/\!\!\sim$ is the
unique one such that the quotient map $\cF\to\cF\!/\!\sim$
defines a morphism
$$
q_\cF:(\cF,\leq)\to(\cF\!/\!\sim,\leq)
$$
of preordered sets, and it is easily seen that the
rule $(\cF,\leq)\mapsto(\cF\!/\!\!\sim,\leq)$ defines
a left adjoint to the inclusion functor
$$
\Poset\to\mathbf{Preorder}.
$$
Moreover, the rule $(\cF,\leq)\mapsto q_\cF$ is a
unit for this adjunction.
\end{example}

The following observations are borrowed from, and
are further developed in \cite[\S I.6]{Gray}.

\begin{remark}\label{rem_adjoint-transf}
(i)\ \
Consider two adjoint pairs $(G_1,F_1)$ and $(G_2,F_2)$ :
$$
\xymatrix{ \cA \ar@<.5ex>[r]^-{F_1}
& \cB \ar@<.5ex>[r]^-{F_2} \ar@<.5ex>[l]^-{G_1} &
\cC \ar@<.5ex>[l]^-{G_2}
}$$
and suppose that for $i=1,2$ we are given adjunctions
$\vartheta_{i,\bullet\bullet}$ for the pair $(G_i,F_i)$.
Then clearly $(G_1\circ G_2,F_2\circ F_1)$ is an adjoint
pair, and we get an induced adjunction for this pair,
by the composition :
$$
\Hom_\cA(G_1G_2C,A)\xrightarrow{\ \vartheta_{1,A,G_2C}\ }
\Hom_\cB(G_2C,F_1A)\xrightarrow{\ \vartheta_{2,F_1A,C}\ }
\Hom_\cB(C,F_2F_1A)
$$
for every $A\in\Ob(\cA)$ and $C\in\Ob(\cC)$. We denote
this adjunction by
$$
(\vartheta_2\circ\vartheta_1)_{\bullet\bullet}
$$
and we call it the {\em composition} of the adjunctions
$\vartheta_1$ and $\vartheta_2$. If $(\eta_i,\eps_i)$ are
the units and counits of $\vartheta_{i,\bullet\bullet}$
(for $i=1,2$), then the unit and counit of
$(\vartheta_2\circ\vartheta_1)_{\bullet\bullet}$
are respectively :
$$
(F_2*\eta_1*G_2)\odot\eta_2
\qquad\text{and}\qquad
\eps_1\odot(G_1*\eps_2*F_1).
$$
Moreover, suppose that $\xymatrix{ \cC \ar@<.5ex>[r]^-{F_3}
& \cD \ar@<.5ex>[l]^-{G_3}}$ is another adjoint pair of
functors, and $\theta_{3,\bullet\bullet}$ an adjunction for
this pair; with this notation, it is also then clear that
$$
(\theta_3\circ(\theta_2\circ\theta_1))_{\bullet\bullet}=
(\theta_3\circ(\theta_2\circ\theta_1))_{\bullet\bullet}.
$$

(ii)\ \
Suppose that we have two pairs of adjoint functors
and two natural transformations
$$
\xymatrix{ \cA \ar@<.5ex>[r]^-F
& \cB \ar@<.5ex>[l]^-G &
\cA \ar@<.5ex>[r]^-{F'} & \cB \ar@<.5ex>[l]^-{G'}
}\qquad
\tau:F\Rightarrow F'
\qquad
\mu:G'\Rightarrow G
$$
and let us fix units and counits $(\eta,\eps)$ (resp.
$(\eta',\eps')$) for the adjoint pair $(G,F)$ (resp.
$(G',F')$). Then we obtain {\em adjoint transformations}
$$
\tau^\dagger:G'\Rightarrow G
\qquad
\mu^\dagger:F\Rightarrow F'
$$
given by the compositions :
$$
\begin{aligned}
G'B & \xrightarrow{\ G'(\eta_B) }G'FGB\xrightarrow{\ G'(\tau_{GB}) }
G'F'GB\xrightarrow{\ \eps'_{GB} }GB
\qquad & &
\text{for every $B\in\Ob(\cB)$} \\
FA & \xrightarrow{\ \eta'_{FA}\ }F'G'FA
\xrightarrow{\ F'(\mu_{FA})\ }
F'GFA\xrightarrow{\ F'(\eps_A)\ }F'A
\qquad & &
\text{for every $A\in\Ob(\cA)$}.
\end{aligned}
$$
We claim that $(\tau^\dagger)^\dagger=\tau$ and
$(\mu^\dagger)^\dagger=\mu$. Indeed, let $\vartheta_{\bullet\bullet}$
(resp. $\vartheta'_{\bullet\bullet}$) be the adjunctions
corresponding to $(\eta,\eps)$ (resp. to $(\eta',\eps')$),
and notice that
$$
\tau^\dagger_B=\vartheta^{\prime-1}_{GB,FGB}(\tau_{GB})\circ G'(\eta_B)
\qquad
\mu^\dagger_A=F'(\eps_A)\circ\vartheta'_{GFA,FA}(\mu_{FA})
$$
so we may compute :
$$
\begin{aligned}
(\tau^\dagger)^\dagger_A
=\,& F'(\eps_A)\circ\vartheta'_{GFA,FA}
(\vartheta^{\prime-1}_{GFA,FGFA}(\tau_{GFA})\circ G'(\eta_{FA})) \\
=\,& \vartheta'_{A,FA}(\eps_A\circ\vartheta^{\prime-1}_{GFA,FGFA}
(\tau_{GFA})\circ G'(\eta_{FA})) \\
=\,& \vartheta'_{A,FA}
(\vartheta^{\prime-1}_{A,FGFA}(F'(\eps_A)\circ\tau_{GFA})
\circ G'(\eta_{FA})) \\
=\,& \vartheta'_{A,FA}
(\vartheta^{\prime-1}_{A,FGFA}(\tau_A\circ F(\eps_A))\circ G'(\eta_{FA})) \\
=\,& \vartheta'_{A,FA}
(\vartheta^{\prime-1}_{A,FA}(\tau_A)\circ G'F(\eps_A)\circ G'(\eta_{FA})) \\
=\,& \vartheta'_{A,FA}(\vartheta^{\prime-1}_{A,FA}(\tau_A)) \\
=\,& \tau_A
\end{aligned}
$$
where the second, third and fifth identities follow from the
naturality of $\vartheta'_{\bullet\bullet}$, the fourth from
the naturality of $\tau$, and the sixth from the triangular
identities of \eqref{subsec_adj-pair}. We leave to the reader
the similar calculation which gives the second identity.
Hence the rule
$$
\tau\mapsto\tau^\dagger:=(\tau,\vartheta,\vartheta')^\dagger
$$
establishes a natural bijection from the set of natural
transformations $F\Rightarrow F'$, to the set of natural
transformations $G'\Rightarrow G$.
Notice that this correspondence depends not only on $(G,F)$
and $(G',F')$, but also on $(\eta,\eps)$ and $(\eta',\eps')$.
Sometimes we denote this adjoint transformation also by
$(\tau,\eta,\eta')^\dagger$.

(iii)\ \
Moreover, using the triangular identities of
\eqref{subsec_adj-pair}, it is easily seen that the diagrams :
$$
\xymatrix@C+10pt{ G'\circ F \ar@{=>}[r]^-{G'*\tau}
\ar@{=>}[d]_{\tau^\dagger*F} & 
G'\circ F' \ar@{=>}[d]^{\eps'} &
\one_\cB \ar@{=>}[r]^-\eta \ar@{=>}[d]_{\eta'} &
F\circ G \ar@{=>}[d]^{\tau*G} \\
G\circ F \ar@{=>}[r]^-\eps & \one_\cA &
F'\circ G' \ar@{=>}[r]^-{F'*\tau^\dagger} & F'\circ G
}$$
commute. Also, $(\tau,\vartheta,\vartheta')^\dagger$ is
characterized as the unique natural transformation
$G'\Rightarrow G$ such that the following diagram
commutes for every $A\in\Ob(\cA)$ and $B\in\Ob(\cB)$ :
\set\begin{equation}\label{eq_dag-adjoints}
{\diagram \Hom_\cA(GB,A)
\ar[rr]^-{\vartheta_{A,B}} \ar[d]_{\Hom_\cA(\tau^\dagger_B,A)} & &
\Hom_\cB(B,FA) \ar[d]^{\Hom_\cB(B,\tau_A)} \\
\Hom_\cA(G'B,A) \ar[rr]^-{\vartheta'_{A,B}} & & \Hom_\cB(B,F'A)
\enddiagram}
\end{equation}
Indeed, letting $A:=GB$ and recalling that
$\vartheta_{GB,B}(\one_{GB})=\eta_B$, we see easily that the
commutativity of \eqref{eq_dag-adjoints} determines uniquely
$\tau^\dagger$ (details left to the reader). Conversely, if
$\tau^\dagger$ is defined as in (i), we may compute, for every
morphism $f:GB\to A$ in $\cA$ :
$$
\begin{aligned}
\vartheta'_{A,B}(f\circ\tau^\dagger_B)
=\,& \vartheta'_{A,B}
(f\circ\vartheta^{\prime-1}_{GB,FGB}(\tau_{GB})\circ G'(\eta_B)) \\
=\,& \vartheta'_{A,FGB}(f\circ\vartheta^{\prime-1}_{GB,FGB}(\tau_{GB}))
\circ\eta_B \\
=\,& \vartheta'_{A,FGB}(\vartheta^{\prime-1}_{A,FGB}(F'f\circ\tau_{GB}))
\circ\eta_B \\
=\,& F'f\circ\tau_{GB}\circ\eta_B \\
=\,& \tau_A\circ Ff\circ\eta_B \\
=\,& \tau_A\circ\vartheta_{A,B}(f).
\end{aligned}
$$

(iv)\ \
Furthermore, suppose we have a third pair of adjoint functors
$$
\xymatrix{ \cA \ar@<.5ex>[r]^-{F''}
& \cB \ar@<.5ex>[l]^-{G''}
}
\qquad
\text{and a natural transformation}
\qquad
\omega:F'\Rightarrow F''
$$
and let us fix an adjunction $\vartheta''_{\bullet\bullet}$
for the pair $(G'',F'')$. Then we have :
$$
(\omega\odot\tau,\vartheta,\vartheta'')^\dagger=
(\tau,\vartheta,\vartheta')^\dagger\odot
(\omega,\vartheta',\vartheta'')^\dagger.
$$
Indeed, this identity follows easily from the
characterization of $\tau^\dagger$, $\omega^\dagger$
and $(\omega\circ\tau)^\dagger$ given in (iii) (details
left to the reader).

(v)\ \
Lastly, in the situation of (i), suppose moreover that
we have two other adjoint pairs
$$
\xymatrix{
\cA \ar@<.5ex>[r]^-{F'_1}
& \cB \ar@<.5ex>[l]^-{G'_1} \ar@<.5ex>[r]^-{F'_2} &
\cC \ar@<.5ex>[l]^-{G'_2}}
\qquad
\text{and natural transformations}
\qquad
\tau_1:F_1\Rightarrow F'_1
\qquad
\tau_2:F_2\Rightarrow F_2'
$$
and for $i=1,2$, let us fix an adjunction
$\vartheta'_{i,\bullet\bullet}$ for $(G'_i,F'_i)$.
Then we get as in (ii) the natural transformations
$\tau_1^\dagger:G'_1\Rightarrow G_1$ and
$\tau_2^\dagger:G'_2\Rightarrow G_2$, and we have
the identity
$$
(\tau_2*\tau_1,\vartheta_2\circ\vartheta_1,
\vartheta'_2\circ\vartheta'_1)^\dagger=
(\tau_1,\vartheta_1,\vartheta'_1)^\dagger*
(\tau_2,\vartheta_2,\vartheta'_2)^\dagger.
$$
Indeed, taking into account (iii) we get the commutative
diagram
$$
\xymatrix{
\Hom_\cA(G_1G_2C,A) \ar[rr]^-{\vartheta_{1,A,G_2C}}
\ar[d]_{\Hom_\cA(\tau^\dagger_{1,G_2C},A)} & &
\Hom_\cB(G_2C,F_1A) \ar[rr]^-{\vartheta_{2,F_1A,C}}
\ar[d]_{\Hom_\cB(G_2C,\tau_{1,C})} & &
\Hom_\cC(C,F_2F_1A) \ar@<-2ex>[d]^{\Hom_\cC(C,F_2(\tau_{1,C}))} \\
\Hom_\cA(G'_1G_2C,A) \ar[rr]^-{\vartheta'_{1,A,G_2C}}
\ar[d]_{\Hom_\cA(G'_1(\tau^\dagger_{1,C}),A)} & &
\Hom_\cB(G_2C,F'_1A) \ar[rr]^-{\vartheta_{2,F'_1A,C}}
\ar[d]_{\Hom_\cB(\tau^\dagger_{2,C},F'_1A)} & &
\Hom_\cC(C,F_2F'_1A) \ar@<-2ex>[d]^{\Hom_\cC(C,\tau_{2,F'_1C})} \\
\Hom_\cA(G'_1G'_2C,A) \ar[rr]^-{\vartheta'_{1,A,G'_2C}} & &
\Hom_\cB(G'_2C,F'_1A) \ar[rr]^-{\vartheta'_{2,F'_1A,C}} & &
\Hom_\cC(C,F'_2F'_1A)
}$$
for every $A\in\Ob(\cA)$ and $C\in\Ob(\cC)$. The
sought identity follows after invoking again (iii).

(vi)\ \
Especially, taking into account the triangular identities
of \eqref{subsec_adj-pair}, it is easily seen that :
$$
\begin{aligned}
(\one_{F_2},\vartheta_2,\vartheta_2)^\dagger&\,=\one_{G_2} \\
(F_2*\tau_1,\vartheta_2\circ\vartheta_1,
\vartheta_2\circ\vartheta'_1)^\dagger&\,=
(\tau_1,\vartheta_1,\vartheta'_1)^\dagger*G_2 \\
(\tau_2*F_1,\vartheta_2\circ\vartheta_1,
\vartheta'_2\circ\vartheta_1)^\dagger&\,=
G_1*(\tau_2,\vartheta_2,\vartheta'_2)^\dagger.
\end{aligned}
$$
\end{remark}

\begin{remark}\label{rem_opposite-Fun}
(i)\ \
For any two categories $\cC,\cD$ we have a
natural isomorphism of categories
$$
\bFun(\cC,\cD)^o\isom\bFun(\cC^o,\cD^o)
$$
that assigns to any functor $H:\cC\to\cD$ the opposite
functor $H^o:\cC^o\to\cD^o$, and to any natural transformation
$\tau:H\Rightarrow K$ the {\em opposite transformation}
$\tau^o:K^o\Rightarrow H^o$ such that $\tau^o_{C^o}:=(\tau_C)^o$
for every $C\in\Ob(\cC)$. Also, in the situation of
\eqref{subsec_Godem-prod}, notice the identities
$$
(\beta\odot\alpha)^o=\alpha^o\odot\beta^o
\qquad\text{and}\qquad
(\alpha'*\alpha)^o=\alpha'{}^o*\alpha^o.
$$

(ii)\ \
Let $\cB$ and $\cE$ be two other categories, $f:\cB\to\cC$
and $g:\cD\to\cE$ two functors; we get an induced functor
$$
\bFun(f,g):\bFun(\cC,\cD)\to\bFun(\cB,\cE)
\qquad
H\mapsto g\circ H\circ f
\qquad
(\alpha:H\Rightarrow K)\mapsto g*\alpha*f.
$$
Likewise, if $f':\cB\to\cC$ and $g':\cD\to\cE$ are any
two other functors, every pair of natural transformations
$\gamma:f\Rightarrow f'$ and $\delta:g\Rightarrow g'$
induces a transformation
$$
\bFun(\gamma,\delta):\bFun(f,g)\Rightarrow\bFun(f',g')
\qquad
H\mapsto \delta*H*\gamma.
$$
We shall usually write $\bFun(f,\cD)$ (resp. $\bFun(\cC,g)$)
instead of $\bFun(f,\one_\cD)$ (resp. of $\bFun(\one_\cC,g)$).
Furthermore, in the situation of \eqref{subsec_Godem-prod},
notice the identities
$$
\begin{aligned}
\bFun(\beta,\cD)\odot\bFun(\alpha,\cD)
&=\bFun(\beta\odot\alpha,\cD) \\
\bFun(\alpha,\cD)*\bFun(\alpha',\cD)
&=\bFun(\alpha'*\alpha,\cD).
\end{aligned}
$$

(iii)\ \
In the situation of (ii), suppose that the functor
$f$ admits a left adjoint $g:\cC\to\cB$. Then
$f^*:=\bFun(f,\cD)$ is left adjoint to $g^*:=\bFun(g,\cD)$.
More precisely, let $\eta$ be a unit and $\eps$
a counit for the adjoint pair $(g,f)$. From the
triangular identities \eqref{subsec_adj-pair} for
$(g,f)$, and taking into account (ii), we deduce
commutative diagrams :
$$
\xymatrix{ g^*  \ar@{=>}[rrr]^-{\bFun(\eta,\cD)*g^*}
\ar@{=}[rrrd] & & &
g^*f^*g^* \ar@{=>}[d]^{g^**\bFun(\eps,\cD)} & &
f^* \ar@{=>}[rrr]^-{f^**\bFun(\eta,\cD)}
\ar@{=}[rrrd] & & &
f^*g^*f^* \ar@{=>}[d]^{\bFun(\eps,\cD)*f^*} \\
& & & g^* & & & & & f^*
}$$
which, in light of proposition \ref{prop_triangular-identities}(i),
says that $\bFun(\eta,\cD)$ is a unit and $\bFun(\eps,\cD)$
a counit for the adjoint pair $(f^*,g^*)$.

(iv)\ \
Let $F:\cC\to\cD$ be a functor, and $G:\cD\to\cC$ a
left adjoint to $F$. Then $F^o$ is left adjoint to $G^o$.
More precisely, let $\eta$ be a unit and $\eps$ a counit
for the adjoint pair $(G,F)$; then it follows easily from
(i) and proposition \ref{prop_triangular-identities}(i) that
$\eps^o$ is a unit and $\eta^o$ is a counit for the adjoint
pair $(F^o,G^o)$ : details left to the reader.
\end{remark}

\begin{proposition}\label{prop_fullfaith-adjts}
Let $F:\cA\to\cB$ be a functor.
\begin{enumerate}
\item
The following conditions are equivalent :
\begin{enumerate}
\item
$F$ is fully faithful and has a fully faithful left adjoint.
\item
There exist a functor $G:\cB\to\cA$ and isomorphisms
of functors
$$
G\circ F\isom\one_\cA
\qquad
\one_\cB\isom F\circ G.
$$
\item
$F$ is an equivalence.
\end{enumerate}
\item
Suppose that $F$ admits a left adjoint $G:\cB\to\cA$, and let
$\eta:\one_\cB\Rightarrow F\circ G$ and
$\eps:G\circ F\Rightarrow\one_\cA$ be a unit and respectively
a counit for the adjoint pair $(G,F)$. Then $F$ (resp. $G$) is
faithful if and only if $\eps_X$ is an epimorphism for every
$X\in\Ob(\cA)$ (resp. $\eta_Y$ is a monomorphism for every
$Y\in\Ob(\cB)$).
\item
In the situation of {\em (ii)}, the following conditions are
equivalent :
\begin{enumerate}
\item
$F$ (resp. $G$) is fully faithful.
\item
The counit $\eps$ (resp. the unit $\eta$) is an isomorphism
of functors.
\item
There exists an isomorphism of functors $\eps':G\circ F\isom\one_\cA$
(resp. $\eta':\one_\cB\isom F\circ G$).
\end{enumerate}
Moreover, if {\em(c)} holds, there exists a unique adjunction
$\theta'$ for the pair $(G,F)$ whose counit is $\eps'$ (resp.
whose unit is $\eta'$).
\item
Suppose that $F$ admits both a left adjoint $G:\cB\to\cA$
and a right adjoint $H:\cB\to\cA$. Then $G$ is fully faithful
if and only if $H$ is fully faithful.
\end{enumerate}
\end{proposition}
\begin{proof}(ii): In view of remark \ref{rem_opposite-Fun}(iv),
it suffices to consider the assertion relative to $\eps_\bullet$.
Thus, suppose that $Ff_1=Ff_2$ for two morphisms $f_1,f_2:X\to X'$.
The naturality of $\eps_\bullet$ yields the identity :
$\eps_{X'}\circ GFf_i=f_i\circ\eps_X$ for $i=1,2$, and if $\eps_X$
is an epimorphism, we deduce $f_1=f_2$. Conversely, suppose $F$
is faithful and $f_1\circ\eps_X=f_2\circ\eps_X$; from the
triangular identities we get :
$$
Ff_i=F(\eps_{X'})\circ FGF(f_i)\circ\eta_{FX}=
F(f_i\circ\eps_X)\circ\eta_{FX}
$$
so $Ff_1=Ff_2$ and therefore $f_1=f_2$ which shows that
$\eps_X$ is an epimorphism.

(iii): Again, remark \ref{rem_opposite-Fun}(iv) reduces to
considering the assertion for $\eps_\bullet$. We check first
that (iii.a)$\Rightarrow$(iii.b) : indeed, if $F$ is fully
faithful, for every $X\in\Ob(\cA)$, there exists a morphism
$\beta_X:X\to GFX$ such that $F\beta_X=\eta_{FX}:FX\to FGFX$.
From the triangular identities of \eqref{subsec_adj-pair}
we deduce that
$F(\eps_X\circ\beta_X)=F(\eps_X)\circ\eta_{FX}=\one_{FX}$,
whence $\eps_X\circ\beta_X=\one_X$, since $F$ is faithful.
Next, let $\theta$ be the unique adjunction for the pair
$(G,F)$ whose unit and counit are $\eta$ and $\eps$; by
inspection of the explicit description of $\theta$ in the
proof of proposition \ref{prop_triangular-identities}(i)
we get
$$
\theta_{GFX,FX}(\beta_X\circ\eps_X)=
F(\beta_X\circ\eps_X)\circ\eta_{FX}=
\eta_{FX}\circ F(\eps_X)\circ\eta_{FX}=\eta_{FX}=
\theta_{GFX,FX}(\one_{GFX})
$$
whence $\beta_X\circ\eps_X=\one_{GFX}$. Obviously
(iii.b)$\Rightarrow$(iii.c). Lastly, suppose that (iii.c)
holds; for every $X,Y\in\Ob(\cA)$ we deduce a bijection
$$
\xi_{X,Y}:\Hom_\cA(X,Y)\xrightarrow{\phi_{X,Y}}\Hom_\cA(GFX,Y)
\xrightarrow{\theta_{GX,Y}}\Hom_\cB(FX,FY)
$$
where $\phi_{X,Y}(f):=f\circ\eps'_X$ for every morphism
$f:X\to Y$ in $\cA$. It is easily seen that $\xi_{X,Y}$ is
natural in both $X$ and $Y$; especially, for every
$f\in\Hom_\cA(X,Y)$ we have :
$$
\xi_{Y,Y}(\one_Y)\circ Ff=\xi_{X,Y}(\one_Y\circ f)=
\xi_{X,Y}(f\circ\one_X)=Ff\circ\xi_{X,X}(\one_X).
$$
Letting $X=Y$, and taking $f:X\to X$ with
$\xi_{X,X}(f)=\one_{FX}$, we deduce that $\xi_{X,X}(\one_X)$
is an isomorphism in $\cB$, for every $X\in\Ob(\cA)$.
Since $\xi_{X,Y}(f)=Ff\circ\xi_{X,X}(\one_X)$ for every
$f\in\Hom_\cA(X,Y)$, we conclude that the rule $f\mapsto Ff$
is a bijection $\Hom_\cA(X,Y)\isom\Hom_\cB(FX,FY)$, whence
(iii.a). Lastly, let us check the existence and uniqueness
of the adjunction $\theta'$ whose counit is $\eps'$. To this
aim, let $\omega:\one_\cA\isom\one_\cA$ be the automorphism
such that $\omega\circ\eps=\eps'$, and set
$\eta':=(F*\omega^{-1}*G)\circ\eta$. We compute :
$$
\begin{aligned}
(F*\eps')\odot(\eta'*F)&\,=(F*\omega)\odot(F*\eps)\odot
(F*\omega^{-1}*GF)\odot(\eta*F) \\
&\,=(F*\omega)\odot(F*\omega^{-1})\odot(F*\eps)\odot(\eta*F)=\one_F.
\end{aligned}
$$
Likewise we check that $(\eps'*G)\odot(G*\eta')=\one_G$,
whence the contention, by proposition
\ref{prop_triangular-identities}.

(i): The equivalence (i.a)$\Leftrightarrow$(i.b) follows
immediately from (iii). Moreover, if (i.b) holds, then for
every $Y\in\Ob(\cB)$ we have an isomorphism $Y\isom FGY$,
so $F$ is essentially surjective. Since we have just seen
that (i.b) implies that $F$ is fully faithful, we deduce
as well that (i.b)$\Rightarrow$(i.c).

(i.c)$\Rightarrow$(i.a): We construct a left adjoint
$G:\cB\to\cA$ to $F$ as follows. Since $F$ is essentially
surjective, for every $B\in\Ob(\cB)$ we may find an object
$A\in\Ob(\cA)$ with an isomorphism $\eta_B:B\isom FA$, and
we set $GB:=A$. Next, since $F$ is fully faithful, for every
morphism $g:B\to B'$ in $\cB$ there exists a unique morphism
$f:GB\to GB'$ in $\cA$ which makes commute the diagram
$$
\xymatrix{
B \ar[r]^-g \ar[d]_{\omega_B} & B' \ar[d]^{\omega_{B'}} \\
FGB \ar[r]^-{Ff}  & FGB'
}$$
and we let $Gg:=f$. It is easily seen that these rules yield
a well defined functor $G$ as sought, and $\eta$ is then a
natural isomorphism $\one_\cB\Rightarrow FG$. Moreover from the
full faithfulness of $F$ we deduce that $G$ is fully faithful
as well (details left to the reader). To conclude we remark,
more generally :

\begin{claim}\label{cl_get-adjunction}
Let $F:\cA\to\cB$ and $G:\cB\to\cA$ be two functors,
$\eta:\one_\cB\Rightarrow FG$ an isomorphism of functors,
and suppose that $F$ is fully faithful. Then there exists
a unique adjunction $\vartheta$ for the pair $(G,F)$ whose
unit is $\eta$.
\end{claim}
\begin{pfclaim} From the proof of proposition
\ref{prop_triangular-identities} we know that $\eta$
determines $\vartheta$, by the rule :
$\vartheta_{A,B}(f):=Ff\circ\eta_B$ for every $A\in\Ob(\cA)$,
$B\in\Ob(\cB)$ and every morphism $f:GB\to A$ in $\cA$.
Conversely, our assumptions easily imply that this rule
does yield a natural bijection
$\Hom_\cA(GB,A)\isom\Hom_\cB(B,FA)$, whence the contention
(details left to the reader).
\end{pfclaim}

(iv): Let $\eta:\one_\cB\Rightarrow FG$ (resp.
$\eta':\one_\cA\Rightarrow HF$) be the unit and
$\eps:GF\Rightarrow\one_\cA$ (resp. $\eps':FH\Rightarrow\one_\cB$)
the counit of a given adjunction for the adjoint pair $(G,F)$
(resp. $(F,H)$). By remark \ref{rem_opposite-Fun}(iv) we
may assume that $H$ is fully faithful, and we show that
the same holds for $G$. By (iii), this is the same as
assuming that $\eps'$ is an isomorphism, and we need to
check that the same holds for $\eta$. To this aim, denote
$\gamma:FG\Rightarrow\one_\cB$ the composition
$$
FG\xrightarrow{\ (FG)*\eps'{}^{-1}\ }FGFH
\xrightarrow{\ F*\eps*H\ }FH\xrightarrow{\ \eps'\ }\one_\cB.
$$
We show that $\gamma$ is inverse to $\eta$. Indeed we have
$$
\gamma\odot\eta=
\eps'\odot(F*\eps*H)\odot(\eta*FH)\odot\eps'{}^{-1}=
\eps'\odot\eps'{}^{-1}=\one_\cB
$$
where the first identity holds by the naturality of
$\eta$, and the second follows from the triangular
identites of \eqref{subsec_adj-pair}. Likewise, we have
$$
\begin{aligned}
\eta\odot\gamma=\,&
(\eps'*FG)\odot(FH*\eta)\odot(F*\eps*H)\odot(FG*\eps'{}^{-1}) \\
=\,& (\eps'*FG)\odot(F*\eps*HFG)\odot(FGFH*\eta)\odot(FG*\eps'{}^{-1}) \\
=\,& (\eps'*FG)\odot(F*\eps*HFG)\odot(FG*\eps'{}^{-1}*FG)\odot(FG*\eta) \\
=\,& (\eps'*FG)\odot(F*\eps*HFG)\odot(FGF*\eta'*G)\odot(FG*\eta) \\
=\,& (\eps'*FG)\odot(F*\eta'*G)\odot(F*\eps*G)\odot(FG*\eta) \\
=\,& (\one_F*G)\odot(F*\one_G)=\one_{FG}
\end{aligned}
$$
where the first and third identities follow from the naturality
of $\eps'$, the second and fifth from that of $\eps$, the fourth
and sixth from the triangular identities for the pairs
$(\eta',\eps')$ and $(\eta,\eps)$.
\end{proof}

\begin{definition} Let $F:\cA\to\cB$ be any equivalence
of categories. A {\em quasi-inverse} for $F$ is the datum
of a functor $G:\cB\to\cA$ and an adjunction for the pair
$(G,F)$. Then we also say that $F$ is a {\em quasi-inverse}
for $G$.
\end{definition}

\begin{remark} It follows easily from proposition
\ref{prop_fullfaith-adjts}(i,iii) and claim \ref{cl_get-adjunction},
that a quasi-inverse for an equivalence $F:\cA\to\cB$ is the
same as the datum of a functor $G:\cB\to\cA$ and an isomorphism
of functors $\one_\cB\isom FG$, and moreover any such $G$ is
also an equivalence. Taking into account remark
\ref{rem_opposite-Fun}(iv), we see that a quasi-inverse for
$F$ is also the same as the datum of such a $G$ and an
isomorphism of functors $GF\isom\one_\cA$.
\end{remark}

\sset\subsubsection{Slice categories}\label{subsec_slice-cat}
A standard construction attaches to any $X\in\Ob(\cC)$
a category :
$$
\cC/X
$$
as follows. The objects of $\cC/X$ are all the pairs
$(A,f)$ where $A\in\Ob(\cC)$ and $f:A\to X$ is any
morphism of $\cC$. For any two such objects $(A,f)$, $(B,g)$,
the set $\Hom_{\cC/X}((A,f),(B,g))$ consists of all the
commutative diagrams of morphisms of $\cC$ :
$$
\xymatrix{ A \ar[rr]^-h \ar[dr]_f & & B \ar[dl]^g \\
& X
}$$
with composition of morphisms induced by the composition law
of $\cC$. We denote sometimes such a morphism of $\cC/X$ by
$$
h/X:(A,f)\to(B,g).
$$
An object (resp. a morphism) of $\cC/X$ is also called an
{\em $X$-object} (resp. an {\em $X$-morphism}) of $\cC$.
Dually, one defines
$$
X/\cC:=(\cC^o/X^o)^o
$$
{\em i.e.} the objects of $X/\cC$ are the pairs $(A,f)$
with $A\in\Ob(\cC)$ and $f:X\to A$ any morphism of $\cC$.
We have an obvious faithful {\em source\/} functor
$$
\ss_X:\cC/X\to\cC
\qquad
(A,f)\mapsto A
\qquad
((A,f)\xrightarrow{h/X}(B,g))\mapsto(A\xrightarrow{h}B)
$$
and likewise one obtains a {\em target\/} functor
$$
\st_X:=\ss^o_{X^o}:X/\cC\to\cC.
$$
Moreover, any morphism $f:X\to Y$ in $\cC$ induces functors :
\set\begin{equation}\label{eq_push-for}
\begin{aligned}
f_*:\cC/X\to\cC/Y \quad : & \quad
(A,g:A\to X)\mapsto(A,f_*g:=f\circ g:A\to Y) \\
f^*:Y\!/\cC\to X\!/\cC \quad : & \quad
(B,h:Y\to B)\mapsto(B,f^*h:=h\circ f:X\to B).
\end{aligned}
\end{equation}
Furthermore, given a functor $F:\cC\to\cB$, any $X\in\Ob(\cC)$
induces functors :
\set\begin{equation}\label{eq_restrict-over-X}
\begin{aligned}
F_{|X}:\cC\!/\!X\to\cB\!/\!FX \quad : & \quad
(A,g)\mapsto(FA,Fg) \\
{}_{X|}F:X\!/\cC\to FX\!/\!\cB \quad :& \quad
(B,h)\mapsto(FB,Fh).
\end{aligned}
\end{equation}

\sset\subsubsection{}\label{subsec_fibreovercat}
The categories $\cC/X$ and $X/\cC$ are special cases of the following
more general construction. Let $F:\cA\to\cB$ be any functor. For any
$B\in\Ob(\cB)$, we define
$$
F\cA/B
$$
as the category whose objects are all the pairs $(A,f)$, where
$A\in\Ob(\cA)$ and $f:FA\to B$ is a morphism in $\cB$. The
morphisms $g:(A,f)\to(A',f')$ are the morphisms $g:A\to A'$ in
$\cA$ such that $f'\circ Fg=f$. There are well-defined functors :
$$
F\!/\!B:F\cA/B\to\cB/B \ \ : \ \ (A,f)\mapsto (FA,f)
\qquad \text{and} \qquad
\ss_B:F\cA/B\to\cA \ \ : \ \ (A,f)\mapsto A.
$$
Dually, we define :
$$
B/F\cA:=(F^o\cA^o/B^o)^o
$$
and likewise one has natural functors :
$$
B/F:B/F\cA\to B/\cB
\qquad \text{and} \qquad
\st_B:B/F\cA\to\cA.
$$
Any morphism $g:B'\to B$ induces functors :
$$
\begin{aligned}
g/F\cA:B/F\cA\to B'/F\cA \quad : & \quad (A,f)\mapsto(A,f\circ g) \\
F\cA/g:F\cA/B'\to F\cA/B \quad : & \quad (A,f)\mapsto(A,g\circ f).
\end{aligned}
$$
Obviously, the category $\cC\!/\!X$ (resp. $X\!/\cC$) is the same as
$\one_\cC\cC/\!X$ (resp. $X\!/\one_\cC\cC$).

\sset\subsubsection{}\label{subsec_target-fctr}
In the situation of \eqref{subsec_fibreovercat}, the
categories of the form $F\cA/B$ can be faithfully
embedded in a single category $F\cA/\cB$. The latter
is the category whose objects are all the triples
$(A,B,f)$, where $A\in\Ob(\cA)$, $B\in\Ob(\cB)$ are
any two objects, and $f:FA\to B$ is any morphism of $\cB$.
If $f:FA\to B$ and $f':FA'\to B'$ are any two objects,
the set
$$
\Hom_{F\cA/\cB}((A,B,f),(A',B',f'))
$$
consists of all pairs $(g,g')$ where $g$ is a morphism
in $\cA$ and $g'$ a morphism in $\cB$, that make commute the
diagram :
\set\begin{equation}\label{eq_morph-in-morph}
{\diagram FA \ar[r]^-f \ar[d]_{Fg} & B \ar[d]^{g'} \\
          FA' \ar[r]^-{f'} & B'
\enddiagram}
\end{equation}
with composition of morphisms induced by the composition
laws of $\cA$ and $\cB$, in the obvious way.
There are two natural {\em source\/} and {\em target\/} functors :
$$
\cA\xleftarrow{\ \ss\ }F\cA/\cB\xrightarrow{\ \st\ }\cB
$$
such that $\ss(FA\to B):=A$, $\st(FA\to B):=B$ for any object
$FA\to B$ of $F\cA/\cB$, and $\ss(g,g')=g$, $\st(g,g')=g'$.
Dually, we let
$$
\cB/F\cA:=(F^o\cA^o/\cB^o)^o
$$
and the corresponding source and target functors are switched :
$$
\cA\xleftarrow{\ \st\ }\cB/F\cA\xrightarrow{\ \ss\ }\cB.
$$

\sset\subsubsection{}\label{subsec_Morph-cat}
For the special case of the identity endofunctor $\one_\cC$
of any category $\cC$, we obtain the
{\em category of arrows of\/ $\cC$}
$$
\sMorph(\cC):=\one_\cC\cC/\cC.
$$
So, the set of objects of $\sMorph(\cC)$ is $\rMorph(\cC)$
(notation of \eqref{sec_label-for-morph}) and the morphisms
are the commutative square diagrams in $\cC$. The functor
$\ss_X$ of \eqref{subsec_slice-cat} is the restriction of
$\ss:\sMorph(\cC)\to\cC$ to the subcategory $\cC\!/\!X$,
and likewise for $\st_X$. Likewise, in the situation of
\eqref{subsec_target-fctr}, the target functor on $F\cA/\cB$
and the source functor on $\cB/F\cA$ factor through functors
$$
F\cA/\cB\xrightarrow{\ \sT\ }\sMorph(\cB)\xleftarrow{\ \sS\ }
\cB/F\cA
$$
where $\sT(A,B,f):=(FA,B,f)$ for every object $(A,B,f)$ of
$F\cA/\cB$, and $\sT$ assigns to any morphism
$(g,g'):(A,B,f)\to(A',B',f')$ the commutative square
\eqref{eq_morph-in-morph}, regarded as a morphism
$(FA,B,f)\to(FA',B',f')$ in $\sMorph(\cB)$. Likewise one
describes the functor $\sS$.

Notice also the natural transformation
$$
\xymatrix{ \sMorph(\cC) \rrtwocell<5>^\ss_\st{\ \ \ssm} & & \cC
}$$
where $\ssm(A,B,f):=f$ for every $(A,B,f)\in\rMorph(\cC)$.
Furthermore, every functor $F:\cB\to\cC$ induces a functor
$$
\sMorph(F):\sMorph(\cB)\to\sMorph(\cC)
$$
which maps $(A,B,f)\in\rMorph(\cB)$ to
$(FA,FB,Ff)\in\rMorph(\cC)$ and which sends every commutative
square diagram $D$ in $\cB$ to the commutative square diagram
$FD$ in $\cC$.

Notice that a natural transformation $\alpha$ as in
\eqref{eq_Rightarrow} is equivalent to the datum of
a functor
$$
\tilde\alpha:\cA\to\sMorph(\cB)
\qquad\text{such that}\qquad
\ss\circ\tilde\alpha=F
\qquad\text{and}\qquad
\st\circ\tilde\alpha=G.
$$
Namely, one defines $\tilde\alpha$ by the rule :
$A\mapsto(FA,GA,\alpha_A)$ for every $A\in\Ob(\cA)$,
and for every morphism $f:A\to B$ in $\cA$, one lets
$\tilde\alpha(f)$ be the commutative square diagram
\eqref{eq_natural-transform}.

\sset\subsubsection{}\label{subsec_comma-adjunction}
Let $\cA,\cB$ be two categories, $F:\cA\to\cB$ a functor,
$G:\cB\to\cA$ a left adjoint for $F$, and $\theta$ an
adjunction for the pair $(G,F)$. Then for every $A\in\Ob(\cA)$
the functor $F_{|A}:\cA/A\to\cB/FA$ of \eqref{eq_restrict-over-X}
admits a left adjoint that we denote
$$
G_{|A}:\cB/FA\to\cA/A.
$$
Namely, to every $(f:B\to FA)\in\Ob(\cB/FA)$ we assign
the object $(\theta^{-1}_{AB}(f):GB\to A)\in\Ob(\cA/A)$,
and to every morphism $h/FA:(f:B\to FA)\to(f':B'\to FA)$
in $\cB/FA$ we assign the morphism
$Gh/A:\theta^{-1}_{AB}(f)\to\theta^{-1}_{AB'}(f')$ in $\cA/A$.
Indeed, $\theta$ induces an adjunction $\theta_{|A}$ for
the pair $(G_{|A},F_{|A})$ : to every $(f:B\to FA)\in\Ob(\cB/FA)$
and $(g:A'\to A)\in\Ob(\cA/A)$ we assign the bijection
$$
(\theta_{|A})_{g,f}:\Hom_{\cA/A}(\theta^{-1}_{AB}(f),g)\isom
\Hom_{\cB/FA}(f,Fg)
\qquad
h/A\mapsto\theta(h)/FA.
$$
Dually, since $F^o$ is left adjoint to $G^o$ (remark
\ref{rem_opposite-Fun}(iv)), we see that for every $B\in\Ob(\cB)$
the functor ${}_{B|}G:B/\cB\to GB/\cA$  of \eqref{eq_restrict-over-X}
admits the right adjoint
$$
{}_{B|}F:GB/\cA\to B/\cB
\qquad
(GB\xrightarrow{f}A)\mapsto(B\xrightarrow{\theta_{AB}(f)}FA)
\qquad
GB/h\mapsto B/Fh.
$$
The detailed verifications shall be left to the reader.
See example \ref{ex_comma-adjunction} for a related result.

\sset\subsubsection{}\label{subsec_comma-isoms-from-adj}
In the situation of \eqref{subsec_comma-adjunction} we get
furthermore for every $A\in\Ob(\cA)$ an isomorphism of categories :
$$
\theta^A:G\cB/A\isom\cB/FA
\qquad
(f:GB\to A)\mapsto(\theta_{AB}(f):B\to FA)
$$
that assigns to every morphism $g:(f:GB\to A)\to(f':GB'\to A)$
of $G\cB/A$ the morphism $g:(\theta_{AB}(f):B\to FA)\to
(\theta_{AB'}(f'):B'\to FA)$ of $\cB/FA$. Also, for every
morphism $h:A\to A'$ of $\cA$, we get a commutative diagram
of categories :
$$
\xymatrix{ G\cB/A \ar[r]^-{\theta^A} \ar[d]_{G\cB/h} &
\cB/FA \ar[d]^{(Fh)_*} \\
 G\cB/A' \ar[r]^-{\theta^{A'}} & \cB/FA'.
}$$
Clearly, the isomorphisms $\theta^A$ are restrictions of a single
isomorphism of categories :
$$
G\cB/\cA\isom\cB/F\cA.
$$
The detailed verifications shall be again left to the reader.

\subsection{Presheaves and limits}
A very important construction associated with every category
$\cC$ is the category
$$
\cC^\wedge_\sU:=\bFun(\cC^o,\sU\tdu\Set)
$$
whose objects are called the {\em $\sU$-presheaves\/} on $\cC$
(notation of \eqref{subsec_Godem-prod}). We usually drop the
subscript $\sU$, unless we have to deal with more than one
universe. If $\sU'$ is another universe with $\sU\subset\sU'$,
the natural inclusion of categories :
\set\begin{equation}\label{eq_change-universe}
\cC^\wedge_\sU\to\cC^\wedge_{\sU'}
\end{equation}
is fully faithful (verification left to the reader).
For every functor $F:\cB\to\cC$ and every natural
transformation $\alpha:F\Rightarrow G$ we set
\set\begin{equation}\label{eq_pullback-presheaves}
F^\wedge_\sU:=\bFun(F^o,\Set):
\cC^\wedge_\sU\to\cB^\wedge_\sU
\qquad
\alpha^\wedge_\sU:=\bFun(\alpha^o,\Set):
G^\wedge_\sU\Rightarrow F^\wedge_\sU
\end{equation}
(notation of remark \ref{rem_opposite-Fun}(i,ii)).
Again, we shall drop the subscript and write simply
$F^\wedge$ and $\alpha^\wedge$, unless the omission of
$\sU$ may be a source of ambiguities. Clearly, for
every pair of functors $F_1:\cA\to\cB$ and $F_2:\cB\to\cC$
we have
$$
(F_2\circ F_1)^\wedge=F_1^\wedge\circ F_2^\wedge.
$$
Moreover, in the situation of \eqref{subsec_Godem-prod},
remark \ref{rem_opposite-Fun}(i,iii) yields the identities
\set\begin{equation}\label{eq_wedge-and-godement}
(\beta\odot\alpha)^\wedge=\alpha^\wedge\odot\beta^\wedge
\qquad\text{and}\qquad
(\alpha'*\alpha)^\wedge=\alpha^\wedge*\alpha'{}^\wedge.
\end{equation}

\sset\subsubsection{}\label{subsec_yoneda}
If $\cC$ has small $\Hom$-sets (see \eqref{subsec_opposing}),
there is a natural functor
$$
h_\cC:\cC\to\cC^\wedge
$$
called the {\em Yoneda embedding}, which assigns to every
$X\in\Ob(\cC)$ the functor
$$
h_X:\cC^o\to\Set
\qquad
Y\mapsto\Hom_\cC(Y,X)
\quad
\text{for every $Y\in\Ob(\cC)$}
$$
and to any morphism $f:X\to X'$ in $\cC$, the natural transformation
$h_f:h_X\Rightarrow h_{X'}$ such that
$$
h_{f,Y}(g):=f\circ g
\qquad
\text{for every $Y\in\Ob(\cC)$ and every $g\in\Hom_\cC(Y,X)$}.
$$

\begin{definition}
Let $\cC$ be a category, and $\sV$ any universe such
that $\cC$ has $\sV$-small $\Hom$-sets. We say that an
object $F$ of $\cC^\wedge_\sV$ is
{\em representable in $\cC$\/}, if there exists an isomorphism
$$
h_X\isom F
\qquad
\text{in $\cC^\wedge_\sV$}
$$
for some $X\in\Ob(\cC)$, in which case we also say that $X$
{\em represents} $F$. From the full faithfulness of
\eqref{eq_change-universe}, it follows that the representability
of a presheaf is independent of the universe $\sV$.
\end{definition}

\begin{proposition}[Yoneda's lemma]\label{prop_yoneda}
With the notation of \eqref{subsec_yoneda}, we have :
\begin{enumerate}
\item
The functor $h_\cC$ is fully faithful.
\item
Moreover, for every $F\in\Ob(\cC^\wedge)$ and every $X\in\Ob(\cC)$
there is a natural bijection
$$
F(X)\isom\Hom_{\cC^\wedge}(h_X,F)
$$
functorial in both $X$ and $F$.
\end{enumerate}
\end{proposition}
\begin{proof} Clearly it suffices to check (ii). However,
the sought bijection is obtained explicitly as follows.
To a given $a\in F(X)$, we assign the natural transformation
$\tau_a:h_X\Rightarrow F$ such that
$$
\tau_{a,Y}(f):=Ff(a)
\qquad
\text{for every $Y\in\Ob(\cC)$ and every $f\in h_X(Y)$}.
$$
Conversely, to a given natural transformation
$\tau:h_X\Rightarrow F$ we assign $\tau_X(\one_X)\in F(X)$.
It is easily seen that these rules establish mutually inverse
bijections. The functoriality in $F$ is immediate, and the
functoriality in the argument $X$ amounts to the commutativity
of the diagram
$$
\xymatrix{
\Hom_{\cC^\wedge}(h_{X'},F) \ar[r]^-\sim \ar[d]_{\Hom_{\cC^\wedge}(h_\phi,F)} &
F(X') \ar[d]^{F(\phi)} \\
\Hom_{\cC^\wedge}(h_X,F) \ar[r]^-\sim & F(X)
}$$
for every morphism $\phi:X\to X'$ in $\cC$ : the verification
shall also be left to the reader.
\end{proof}

\begin{example}\label{ex_presheaves-on-sets}
(i)\ \
Every representable presheaf on the category $\Set$ admits
a natural representing object. Indeed, let $\bone$ denote
any set with one element (a canonical choice for $\bone$ is
the set $\{\emptyset\}$); then for every representable
$F\in\Ob(\Set^\wedge)$ we get a natural isomorphism
$$
F\isom h_{F(\bone)}
$$
as follows. For any set $S$, we have a natural bijection
$S\isom\Hom_\Set(\bone,S)$ that assigns to every $s\in S$
the unique map $\bone_s:\bone\to S$ whose image is $\{s\}$.
We deduce a map
$$
S\times F(S)\to F(\bone)
\qquad
(s,a)\mapsto F(\bone_s)(a)
$$
which is the same as a map $F(S)\to\Hom_\Set(S,F(\bone))$
that realizes the sought isomorphism.

(ii)\ \
We may apply Yoneda's lemma to prove the uniqueness
of the (left or right) adjoint for a given functor.
Indeed, let $G:\cB\to\cA$ be a functor between any
two categories, suppose that $F$ and $F'$ are both
right adjoint to $G$, and fix adjunctions
$\vartheta_{\bullet\bullet}$ and $\vartheta'_{\bullet\bullet}$
for $F$ and respectively $F'$. Choose also a universe
$\sU$ such that both $\cA$ and $\cB$ are $\sU$-small.
Then $\vartheta$ and $\vartheta'$ can be regarded as
isomorphisms from the functor $G^\wedge\circ h_\cA$ to
$h_\cB\circ F$ and respectively to $h_\cB\circ F'$. In this
situation, proposition \ref{prop_yoneda}(i) implies that
there exists a unique isomorphism $\omega:F\isom F'$
which makes commute the diagram
$$
\xymatrix{
& G^\wedge\circ h _\cA\ar[ld]_\vartheta \ar[rd]^{\vartheta'} \\
h_\cB\circ F \ar[rr]^-{h_\cB*\omega} & & h_\cB\circ F'.
}$$
Taking into account remark \ref{rem_opposite-Fun}(iv),
a dual argument yields a corresponding uniqueness assertion
for left adjoints.
\end{example}

\begin{remark}\label{rem_represent-morph}
(i)\ \
Let $\cC$ be any category, and $F_1,F_2$ two representable
presheaves on $\cC$, and pick $X_1,X_2\in\Ob(\cC)$ with
isomorphisms $\omega_i:h_{X_i}\isom F_i$ ($i=1,2$). It
follows from proposition \ref{prop_yoneda}(i) that for every
morphism $f:F_1\to F_2$ in $\cC^\wedge$ there exists a unique
morphism $\phi:X_1\to X_2$ in $\cC$ that makes commute the
diagram
$$
\xymatrix{ h_{X_1} \ar[r]^-{h_\phi} \ar[d] & h_{X_2} \ar[d] \\
F_1 \ar[r]^-f & F_2.
}$$
In this situation, we say that $\phi$ {\em represents} the
morphism $f$.

(ii)\ \
Likewise, suppose that $G:\cB\to\cC^\wedge$ is any functor.
We say that $G$ is {\em representable} if there exists a
functor $\gamma:\cB\to\cC$ with an isomorphism
$G\isom h_\cC\circ\gamma$. It follows easily from (i) that
$G$ is representable if and only if $G(B)$ is representable
in $\cC$ for every $B\in\Ob(\cB)$. Moreover, any two
representatives for $G$ are isomorphic in $\bFun(\cB,\cC)$.
Furthermore, if $G,G':\cB\to\cC$ are any two representable
functors, and $\gamma,\gamma':\cB\to\cC$ two corresponding
representing functors, then any natural transformation
$t:G\Rightarrow G'$ is {\em represented} uniquely by a
natural transformation $\tau:\gamma\Rightarrow\gamma'$,
by which we mean that $t=h_\cC*\tau$.
\end{remark}

\sset\subsubsection{Limits and colimits}\label{subsec_wishful}
We wish to explain some standard constructions of presheaves
that are used pervasively throughout this work.
Namely, let $I$ be a small category, $\cC$ a category and $X$
any object of $\cC$. We denote by $c_X:I\to\cC$ the
{\em constant functor\/} associated with $X$ :
$$
c_X(i):=X \quad \text{for every $i\in\Ob(I)$} \qquad
c_X(\phi):=\one_X \quad \text{for every $\phi\in\rMorph(I)$}.
$$
Any morphism $f:X'\to X$ induces a natural transformation
$$
c_f:c_{X'}\Rightarrow c_X \qquad \text{by the rule \ \ :\ \
$(c_f)_i:=f$ for every $i\in I$}.
$$
If $F:I\to\cC$ is any functor, a
{\em cone of vertex $X$ and basis $F$} is any natural
transformation $c_X\Rightarrow F$. Dually, a
{\em cocone with vertex $X$ and basis $F$} is a natural
transformation $F\Rightarrow c_X$.

\begin{definition}\label{def_limits}
With the notation of \eqref{subsec_wishful}, let $\sV$ be
any universe with $\sU\subset\sV$, such that $\cC$ has
$\sV$-small $\Hom$-sets, and let $F:I\to\cC$ be any functor.
\begin{enumerate}
\item
The {\em limit\/} of $F$ is the $\sV$-presheaf on $\cC$ denoted
$$
\lim_I F:\cC^o\to\sV\tdu\Set
$$
that assigns to every $X\in\Ob(\cC)$, the set $\lim_I F(X)$
of all cones $c_X\Rightarrow F$. Any morphism $f:X'\to X$ of
$\cC$ induces the map
$$
\lim_I F(f):\lim_I F(X)\to\lim_I F(X')
\qquad
\tau\mapsto\tau\odot c_f
\quad
\text{for every $\tau:c_X\Rightarrow F$}.
$$
Then we also say that $I$ is the {\em indexing category}
for the limit of $F$.
\item
Dually, the {\em colimit\/} of $F$ is the $\sV$-presheaf on $\cC^o$
$$
\colim_IF:=\lim_{I^o}F^o.
$$
\end{enumerate}
\end{definition}

\begin{remark}\label{rem_wishful}
Let $\cC,\cC'$ be any two categories with $\sV$-small $\Hom$-sets
for some universe $\sU$ such that $\sU\subset\sV$, and $I,I'$ two
small categories.

(i)\ \ 
Any diagram of functors
$$
I'\xrightarrow{\ \phi\ }I\xrightarrow{\ F\ }\cC
\xrightarrow{\ H\ }\cC'
$$
induces a natural morphism
$$
\lim_\phi H:
\lim_I F\to H^\wedge(\lim_{I'}\,H\circ F\circ\phi)
\qquad
\text{in $\cC^\wedge_\sV$}
$$
(notation of \eqref{eq_pullback-presheaves}) by ruling that
$$
\lim_\phi H(X)(\tau):=H*\tau*\phi
\qquad
\text{for every $X\in\Ob(\cC)$ and every cone
$\tau:c_X\Rightarrow F$}.
$$
If $\psi:I''\to I'$ and $K:\cC'\to\cC''$ are any two other
functors, with $I''$ also small and $\cC''$ with $\sV$-small
$\Hom$-sets, the resulting diagram in $\cC^\wedge_\sV$ commutes :
\set\begin{equation}\label{eq_respect}
{\diagram
\lim_IF \ar[d]_{\lim_\phi H}
\ar[rr]^-{\lim_{\phi\circ\psi}K\circ H} & &
H^\wedge\circ K^\wedge(\lim_{I''}K\circ H\circ F\circ\phi\circ\psi) \\
H^\wedge(\lim_{I'}H\circ F\circ\phi)
\ar[rru]_-{\ H^\wedge(\lim_\psi K^\wedge)}.
\enddiagram}
\end{equation}
The reader may spell out the corresponding assertions for
colimits.

(ii)\ \
Any diagram of functors and natural transformations
$$
\xymatrix{
I' \rtwocell^{\phi_1}_{\phi_2}{\alpha} & I \rtwocell^{F_1}_{F_2}{g} & \cC
}$$
induces commutative diagrams
$$
\xymatrix{ \lim_IF_1 \ar[rr]^-{\lim_Ig} \ar[d]_{\lim_{\phi_1}\!\one_\cC}
& & \lim_IF_2 \ar[d]_{\lim_{\phi_2}\one_\cC} &
\colim_IF_2 \ar[rr]^-{\colim_Ig} \ar[d]^{\colim_{\phi_2}\one_\cC}
& & \colim_IF_1 \ar[d]^{\colim_{\phi_1}\!\one_\cC}
 \\
\lim_{I'}(F_1\!\circ\!\phi_1) \ar[rr]^{\lim_{I'}(g*\alpha)} & &
\lim_{I'}(F_2\!\circ\!\phi_2) &
\colim_{I'}(F_2\!\circ\!\phi_2) \ar[rr]^{\colim_{I'}(g*\alpha)} & &
\colim_{I'}(F_1\!\circ\!\phi_1)
}$$
where the top (resp. bottom) arrow of the left diagram is
given by the rule :
$$
\tau\mapsto g\odot\tau
\qquad
\text{(\ resp. $\tau\mapsto(g*\alpha)\odot\tau$ \ )}
$$
for every $X\in\Ob(\cC)$ and every cone $\tau:c_X\Rightarrow F_1$
(resp. $\tau:c_X\Rightarrow F_1\circ\phi_1$). The top
(resp. bottom) arrow of the right diagram is given by the rule
$$
\tau\mapsto\tau\odot g
\qquad
\text{(\ resp. $\tau\mapsto\tau\odot(g*\alpha)$ \ )}
$$
for every $X\in\Ob(\cC)$ and every cocone $\tau:F_2\Rightarrow c_X$
(resp. $\tau:F_2\circ\phi_1\Rightarrow c_X$).

(iii)\ \
Suppose that $\lim_IF$ is representable by $L\in\Ob(\cC)$,
so we have an isomorphism
\set\begin{equation}\label{eq_universal-cone}
\omega:h_L\isom\lim_IF
\qquad
\text{in $\cC^\wedge_\sV$}.
\end{equation}
Notice that the cone $\tau:=\omega_L(\one_L)$
determines $\omega$ : indeed the latter assigns, to every
$X\in\Ob(\cC)$, the bijection :
$$
\omega_X:\Hom_\cC(X,L)\isom\lim_IF(X)
\qquad
f\mapsto \tau\odot c_f.
$$
Conversely, we say that a given $\tau:c_L\Rightarrow F$
is a {\em universal cone}, if the induced map $\omega_X$
is a bijection for every $X\in\Ob(\cC)$, in which case
the limit of $F$ is representable by $L$. Clearly, the
universal property for a cone is independent of the choice
of auxiliary universe $\sV$.

Likewise, if $\colim_IF$ is representable by
$C^o\in\Ob(\cC^o)$, the choice of an isomorphism
\set\begin{equation}\label{eq_universal-cocone}
h_{C^o}\isom\colim_IF
\qquad
\text{in $(\cC^o)^\wedge_\sV$}
\end{equation}
induces a {\em universal cocone} $\mu:F\Rightarrow c_C$,
which in turns determines \eqref{eq_universal-cocone} by the rule :
$$
\Hom_\cC(C,X)\isom\colim_IF(X^o)
\qquad
f\mapsto c_f\odot\mu.
$$
Moreover, with the notation of (ii), suppose that the limits
(resp. colimits) of $F_1$ and $F_2$ are representable by objects
$L_1$ and $L_2$ of $\cC$ (resp. $C_1^o$ and $C_2^o$ of $\cC^o$),
and let us fix isomorphisms as \eqref{eq_universal-cone}
and \eqref{eq_universal-cocone} for $F_1$ and $F_2$. Then :
\begin{itemize}
\item
The limit of $g$ is represented by a morphism $L_1\to L_2$
in $\cC$.
\item
The colimit of $g$ is represented by a morphism
$C_2^o\to C_1^o$ in $\cC^o$, {\em i.e.} by a morphism
$C_1\to C_2$ in $\cC$ (see remark \ref{rem_represent-morph}(ii)).
\end{itemize}

(iv)\ \
Suppose that the universe $\sV$ is sufficiently large,
so that both $\cC$ and $\cC^\wedge_\sU$ have $\sV$-small
$\Hom$-sets. Every $i\in I$ induces a morphism of presheaves
$$
t_i:\lim_IF\to h_{Fi}
$$
that assigns to every $X\in\Ob(\cC)$ and every cone
$\tau:c_X\Rightarrow F$ the morphism $\tau_i:X\to Fi$,
regarded as an element of $h_{Fi}(X)$. We obtain in
this way a natural cone
\set\begin{equation}\label{eq_tautological-cone}
c_{\lim_IF}\Rightarrow h_\cC\circ F
\qquad
i\mapsto t_i
\end{equation}
where $h_\cC:\cC\to\cC^\wedge_\sV$ is the Yoneda embedding.
We call \eqref{eq_tautological-cone} the {\em tautological cone}
associated with $F$. We claim that the tautological cone
is universal. Indeed, say that $G$ is any presheaf on $\cC$,
and $\tau:c_G\Rightarrow h_\cC\circ F$ any cone. We define a
morphism $f_\tau:G\to\lim_IF$ as follows. To every $i\in\Ob(I)$,
the cone $\tau$ attaches a morphism of presheaves
$\tau_i:G\to h_{Fi}$, which is the datum, for every
$X\in\Ob(\cC)$ of a map $\tau_i^X:G(X)\to\Hom_\cC(X,Fi)$,
and for fixed $s\in G(X)$, the system
$$
\tau^X_\bullet(s):=(\tau^X_i(s):X\to Fi~|~X\in\Ob(\cC))
$$
is a cone $c_X\Rightarrow F$. Then we let
$f_\tau(X):G(X)\to\lim_IF(X)$ be the map such that
$s\mapsto\tau_\bullet^X(s)$ for every $X\in\Ob(\cC)$ and
every $s\in G(X)$. The rule $\tau\mapsto f_\tau$ yields an
inverse for the map
$$
\omega_G:\Hom_{\cC^\wedge}(G,\lim_IF)\to\lim_I h_\cC\circ F
$$
associated, as in (iii), with the cone
\eqref{eq_tautological-cone}, whence the claim : details
left to the reader.

Dually, with $F$ we may also associate a
{\em tautological cocone}
$$
h_{\cC^o}\circ F^o\Rightarrow h_{\colim_IF}
\qquad
\text{in $(\cC^o)^\wedge_\sV$}
$$
which is a universal cocone.
\end{remark}

\begin{example}\label{ex_equalizers}
Let $\sV$ be a universe containing $\sU$, and $\cC$ any category
with $\sV$-small $\Hom$-sets.

(i)\ \  
For $i=1,2$, let $f_i:A\to B_i$ be two morphisms in $\cC$; the
{\em push-out\/} or {\em amalgamated sum} of $f_1$ and $f_2$
is the colimit of the functor $F:I\to\cC$, defined as follows.
The set $\Ob(I)$ consists of three objects $s$, $t_1$, $t_2$
and $\rMorph(I)$ consists of two morphisms
$$
t_1\xleftarrow{\ \phi_1\ } s\xrightarrow{\ \phi_2\ } t_2
$$
(in addition to the identity morphisms of the objects of $I$);
the functor is given by the rule : $Fs:=A$, $Ft_i:=B_i$ and
$F\phi_i:=f_i$ (for $i=1,2$).
For any $C\in\Ob(\cC)$, a cocone $F\Rightarrow c_C$ amounts
to a pair of morphisms $f'_i:B_i\to C$ ($i=1,2$) such that
$f'_1\circ f_1=f'_2\circ f_2$. If the cocone is universal
(and thus, $C$ represents the push-out of $f_1$ and $f_2$),
we say that the resulting commutative diagram :
$$
\xymatrix{ A \ar[r]^-{f_1} \ar[d]_{f_2} & B_1 \ar[d]^{f'_1} \\
           B_2 \ar[r]^-{f'_2} & C
}$$
is {\em cocartesian}. Dually one defines the {\em fibre product\/}
or {\em pull-back\/} of two morphisms $g_i:A_i\to B$ ($i=1,2$).
If $D\in\Ob(\cC)$ represents this fibre product, a cone with
vertex $D$ is given by a pair of morphisms $g'_i:D\to A_i$
($i=1,2$) such that $g_1\circ g'_1=g_2\circ g'_2$, and we say
that the diagram :
$$
\xymatrix{ D \ar[r]^-{g'_1} \ar[d]_{g'_2} & A_1 \ar[d]^{g_1} \\
           A_2 \ar[r]^-{g_2} & B
}$$
is {\em cartesian} if this cone is universal. The push-out
of $f_1$ and $f_2$ is a $\sV$-presheaf, and is usually just
called the amalgamated sum of $B_1$ and $B_2$ over $A$,
denoted
$$
B_1\amalg_{(f_1,f_2)}B_2
\qquad\text{or simply}\qquad
B_1\amalg_AB_2
$$
unless the notation gives rise to ambiguities. Likewise one
writes
$$
A_1\times_{(g_1,g_2)}A_2
\qquad\text{or just}\qquad
A_1\times_BA_2
$$
for the fibre product of $g_1$ and $g_2$, which is a $\sV$-presheaf
as well.

(ii)\ \ 
Similarly, consider the category $I'$ with $\Ob(I')=\{s,t\}$,
and with $\rMorph(I')$ consisting of two morphisms
$$
\xymatrix{
s \ar@<.5ex>[r]^-{\psi_1} \ar@<-.5ex>[r]_-{\psi_2} & t
}$$
(in addition to $\one_s$ and $\one_t$).
If $A,B\in\Ob(\cC)$ are any two objects and $f_1,f_2:A\to B$
any two morphisms, we get a functor $F':I'\to\cC$ by the
rule : $s\mapsto A$, $t\mapsto B$ and $\psi_i\mapsto f_i$
for $i=1,2$. Then, the colimit of $F'$ is also called the
{\em coequalizer\/} of $f_1$ and $f_2$, and is sometimes
denoted
$$
\mathrm{Coequal}(f_1,f_2).
$$
Dually, the limit of $F'$ is also called the {\em equalizer\/}
of $f_1$ and $f_2$, sometimes denoted
$$
\Equal(g_1,g_2).
$$

(iii)\ \ 
Let $I\in\sU$ be any small set, and $\underline B:=(B_i~|~i\in I)$
any family of objects of $\cC$. We may regard $I$ as a discrete
category (see example \ref{ex_universe}(ii)), and then the rule
$i\mapsto B_i$ yields a functor $I\to\cC$, whose limit (resp.
colimit) is called the {\em product} (resp. {\em coproduct})
of the family $\underline B$, and is denoted
$$
\prod_{i\in I}B_i
\qquad
\text{(resp.\ \ \  $\coprod_{i\in I}B_i$\ )}.
$$
If $I=\{1,2\}$ is a set with exactly two elements, we also
write $B_1\times B_2$ (resp. $B_1\amalg B_2$) for this limit
(resp. colimit), and we call it sometimes a {\em binary product}
(resp. a {\em binary coproduct}).

(iv)\ \
Let $\cB$ be any category, and suppose that a given
$X\in\Ob(\cB)$ represents the {\em empty product\/}
in $\cB$, {\em i.e.} the product of an empty family
of objects of $\cB$ (this product is a $\sV$-presheaf,
for any suitably large universe $\sV$). This means
precisely that $\Hom_\cB(Y,X)$ consists of exactly
one element, for every $Y\in\Ob(\cB)$. Any such $X$
is called a {\em final object of\/ $\cB$}.

Dually, we say that $X$ is an {\em initial object of\/ $\cB$},
if $\Hom_\cB(X,Y)$ consists of a single element, for every
$Y\in\Ob(\cB)$. Then an object of $\cB$ is initial if and
only if it represents the {\em empty coproduct} in $\cB$.

Moreover, we shall say that an object $X$ of $\cB$
is {\em disconnected}, if there exist $A,B\in\Ob(\cB)$,
neither of which is an initial object of $\cB$, and such
that $X$ represents the coproduct $A\amalg B$ (which, again,
is well defined in any sufficiently large universe $\sV$).
We say that $X$ is {\em connected}, if $X$ is not
disconnected and is not an initial object of $\cB$.

(v)\ \
For instance, the initial object of $\Set$ is the
empty set, and any set $\bone$ with one element is
a final object of $\Set$. More generally, if $\cB$
is any category, the initial object of $\cB^\wedge$
is the presheaf $\emptyset_\cB$ such that
$\emptyset_\cB(X)=\emptyset$ for every $X\in\Ob(\cB)$,
and the presheaf $\bone_\cB$ such that $\bone_\cB(X)=\bone$
for every $X\in\Ob(\cB)$ is a final object.
\end{example}

\sset\subsubsection{}\label{subsec_codiagonal}
Let $\cC$ be a category, $f:A\to B$ a morphism in $\cC$;
as explained in example \ref{ex_equalizers}(i), the pair of
morphisms $B\xleftarrow{\one_B}B\xrightarrow{\one_B}B$ can
be regarded as an element $\tau\in(B\amalg_{(f,f)}B)(B)$, which,
by Yoneda's lemma (proposition \ref{prop_yoneda}(ii)) corresponds
to a morphism
$$
\pi^o_f:h_B\to B\amalg_AB
\qquad
\text{in $(\cC^o)^\wedge$}.
$$
If $B\amalg_AB$ is representable in $\cC$, this is the same
as a morphism
$$
\pi_f:B\amalg_AB\to B
\qquad
\text{in $\cC$}
$$
called the {\em codiagonal of $f$}. Dually, we have a natural
morphism :
$$
\iota_f:h_A\to A\times_BA
\qquad
\text{in $\cC^\wedge$}
$$
which -- in case $A\times_BA$ is representable in $\cC$ --
is the same as a {\em diagonal morphism of $f$}
$$
\iota_f:A\to A\times_BA
\qquad
\text{in $\cC$}.
$$

\begin{proposition}\label{prop_diagonal-monos}
{\em(i)}\ \
With the notation of \eqref{subsec_codiagonal}, the
following conditions are equivalent :
\begin{enumerate}
\alphaenu
\item
$f:A\to B$ is a monomorphism in $\cC$
\item
$\iota_f:h_A\to A\times_BA$ is an isomorphism in $\cC^\wedge$.
\end{enumerate}
In case $A\times_BA$ is representable in $\cC$, denote
by $A\xleftarrow{p_1}A\times_BA\xrightarrow{p_2}A$ the universal
cone of this fibre product. Then {\em(a)} and {\em(b)} are
moreover equivalent to :
\begin{enumerate}
\alphaenu\addenu\addenu
\item
$\iota_f:A\to A\times_BA$ is an isomorphism in $\cC$.
\item
$\iota_f:A\to A\times_BA$ is an epimorphism in $\cC$.
\item
$p_1=p_2$.
\end{enumerate}

{\em(ii)}\ \
Dually, the following conditions are equivalent :
\begin{enumerate}
\alphaenu
\item
$f:A\to B$ is an epimorphism in $\cC$.
\item
$\pi_f^o:h_B\to B\amalg_AB$ is an isomorphism in $(\cC^o)^\wedge$.
\end{enumerate}
In case $B\amalg_AB$ is representable in $\cC$, denote
by $B\xrightarrow{e_1}B\amalg_AB\xleftarrow{e_2}B$ the
universal cocone of this amalgamated sum. Then {\em(a)}
and {\em(b)} are moreover equivalent to :
\begin{enumerate}
\alphaenu\addenu\addenu
\item
$\pi_f:B\amalg_AB\to B$ is an isomorphism in $\cC$.
\item
$\pi_f:B\amalg_AB\to B$ is a monomorphism in $\cC$.
\item
$e_1=e_2$.
\end{enumerate}
\end{proposition}
\begin{proof}(i.a)$\Leftrightarrow$(i.b): For every $X\in\Ob(\cC)$
the map $\iota_{f,X}:h_A(X)\to A\times_BA(X)$ is given by the
rule : $(X\xrightarrow{g}A)\mapsto(A\xleftarrow{g}X\xrightarrow{g}A)$.
Hence, $\iota_f$ is an isomorphism if and only if for
every $X\in\Ob(\cC)$ and every pair of morphisms $g_1,g_2:X\to A$
such that $f\circ g_1=f\circ g_2$, we have $g_1=g_2$. This means
precisely that $f$ is a monomorphism.

(i.b)$\Leftrightarrow$(i.c) is clear, since the Yoneda imbedding
$h_\cC$ is fully faithful.

(i.c)$\Rightarrow$(i.d) is trivial.

(i.d)$\Rightarrow$(i.e): By definition,
$p_1\circ\iota_f=\one_A=p_2\circ\iota_f$. If $\iota_f$ is
an epimorphism, then $p_1=p_2$.

(i.e)$\Rightarrow$(i.c): Set $p:=p_1$; since
$p\circ\iota_f=\one_A$, it suffices to show that
$\iota_f\circ p=\one_{A\times_BA}$. Due to the universality
of the cone $A\xleftarrow{p_1}A\times_BA\xrightarrow{p_2}A$,
we are then reduced to checking that
$p_1\circ\iota_f\circ p=p_1$ and $p_2\circ\iota_f\circ p=p_2$,
which is clear, since by assumption $p=p_2$.

Assertion (ii) follows from (i) by duality.
\end{proof}

Some special kinds of limits occur frequently in
applications; we gather some of them in the following :

\begin{definition}\label{def_filtered-cols}
Let $I$ be a category.

(i)\ \
We say that $I$ is {\em finite\/} if both $\Ob(I)$ and
$\rMorph(I)$ are finite sets.

(ii)\ \
We say that $I$ is {\em connected}, if
$\Ob(I)\neq\emptyset$ and every two objects $i,j$ can
be connected by a finite sequence of morphisms in $I$ :
\set\begin{equation}\label{eq_zig-zag}
i\to k_1\leftarrow k_2\to\cdots\leftarrow k_n\to j.
\end{equation}

(iii)\ \
We say that $I$ is {\em directed}, if for every
$i,j\in\Ob(I)$ there exist $k\in\Ob(\cC)$ and
morphisms $i\to k\leftarrow j$ in $I$. We say that
$I$ is {\em codirected}, if $I^o$ is directed.

(iv)\ \
We say that $I$ is {\em locally directed} if, for every
$i\in\Ob(I)$, the category $i/I$ is directed. We say
that $I$ is {\em locally codirected}, if $I^o$ is
locally directed.

(v)\ \
We say that $I$ is {\em pseudo-filtered\/} if it is locally
directed, and the following {\em coequalizing condition} holds.
For any $i,j\in\Ob(I)$, and any two morphisms $f,g:i\to j$,
there exist $k\in\Ob(I)$ and a morphism $h:j\to k$ such
that $h\circ f=h\circ g$.

(vi)\ \
We say that $I$ is {\em filtered}, if it is pseudo-filtered
and connected. We say that $I$ is {\em cofiltered\/} if
$I^o$ is filtered.

(vii)\ \
Let $F:I\to\cC$ be a functor to any category $\cC$. Then
$\lim_IF$ is a well defined $\sV$-presheaf, for any sufficiently
large universe $\sV$, and we say that {\em the limit of $F$ is
small} (resp. {\em connected}, resp. {\em codirected}, resp.
{\em locally codirected}, resp. {\em cofiltered}, resp.
{\em finite}), if $I$ is small (resp. connected, resp.
codirected, resp. locally codirected, resp. cofiltered,
resp. finite). Dually, one defines {\em small}, {\em connected},
{\em directed}, {\em locally directed}, {\em filtered\/} and
{\em finite\/} colimits.

(viii)\ \ 
We say that a category $\cC$ is {\em complete\/}
(resp. {\em finitely complete})  if every small limit (resp.
finite limit) of $\cC$ is representable in $\cC$.
We say that $\cC$ is {\em cocomplete} (resp.
{\em finitely cocomplete}) if every small colimit (resp.
finite colimit) of $\cC$ is representable in $\cC^o$.
\end{definition}

\begin{remark}\label{rem_cofinal}
Let $I$ be any category.

(i)\ \
If $I$ is connected and locally directed, then $I$ is
directed. Indeed, let $i,j$ be any two objects of $I$;
we have to show that there exist $k\in\Ob(I)$ and
morphisms $a:i\to k$ and $b:j\to k$. However, by assumption
there exists a sequence \eqref{eq_zig-zag}, for some
$n\in\N$ and objects $k_1,\dots,k_n$ of $I$. Then, a
simple induction reduces the assertion to the case where
$n=2$. But since $I$ is locally directed, we may then
find $k\in\Ob(I)$ and morphisms $c:k_1\to k$, $b:j\to k$,
so it suffices to take $a:=c\circ d$, where $d:i\to k_1$
is the morphism appearing in \eqref{eq_zig-zag}.
Dually, if $I$ is connected and locally codirected, then
$I$ is codirected.

It follows that a category $I$ is filtered if and only
if it is directed, satisfies the coequalizing condition
of definition \ref{def_filtered-cols}(v), and its set
of objects is non-empty.
 
(ii)\ \ 
Let $I$ be a small category. It is easily seen that
$I$ is connected if and only if it is a connected
object of the category $\bCat$ (see example
\ref{ex_equalizers}(iv)). In general, there
is a natural decomposition in $\bCat$ :
$$
I\isom\coprod_{s\in S}I_s
$$
where each $I_s$ is a connected category, and $S$ is a
small set. We call $\{I_s~|~s\in S\}$ the {\em set of
connected components} of $I$, and we denote it
$$
\pi_0(I).
$$
Especially, every small pseudo-filtered category is
a coproduct of a (small) family of filtered categories.
This decomposition induces natural isomorphisms in $\cA^\wedge$ :
$$
\colim_IF\isom\coprod_{s\in S}\colim_{I_s} F\circ e_s
\qquad
\lim_I F\isom\prod_{s\in S}\lim_{s\in S} F\circ e_s
$$
for any functor $F:I\to\cA$, where $e_s:I_s\to I$ is the
natural inclusion functor, for every $s\in S$. The details
shall be left to the reader. 
\end{remark}

\begin{proposition}\label{prop_complete-criteria}
For any category\/ $\cC$ we have :
\begin{enumerate}
\item
$\cC$ is complete (resp. finitely complete) if
and only if the following conditions hold :
\begin{enumerate}
\item
The product of every small (resp. finite) family of objects
is representable in $\cC$.
\item
All equalizers are representable in $\cC$.
\end{enumerate}
\item
Dually, $\cC$ is cocomplete (resp. finitely cocomplete)
if and only if the following holds :
\begin{enumerate}
\item
The coproduct of any small (resp. finite) family of
objects is representable in $\cC$.
\item
All coequalizers are representable in $\cC$.
\end{enumerate}
\end{enumerate}
\end{proposition}
\begin{proof} Clearly it suffices to show (i). Now,
conditions (a) and (b) obviously hold if $\cC$ is complete.
Conversely, suppose that the conditions hold, let $I$ be any
small (resp. finite) category, and $F:I\to\cC$ any functor.
We regard $\Ob(I)$ (resp. $\rMorph(I)$) as a discrete subcategory
of $I$ (resp. of $\sMorph(I)$ : see example \ref{ex_universe}(ii)),
and we consider the diagram
$$
\xymatrix{ \rMorph(I)
\ar@<.5ex>[r]^-{\mathrm{s}} \ar@<-.5ex>[r]_-{\mathrm{t}} &
\Ob(I) \ar[r]^-{\mathrm{ob}} & I
}$$
where $\mathrm{s}$ and $\mathrm{t}$ are respectively
the restrictions of the functors $\ss$ and $\st$ from
\eqref{subsec_Morph-cat}, and $\mathrm{ob}$ is the
inclusion functor. We let as well $\mathrm{m}:
\mathrm{ob}\circ\mathrm{s}\Rightarrow\mathrm{ob}\circ\mathrm{t}$
be the restriction of the natural transformation $\ssm$ of
\eqref{subsec_Morph-cat}. By assumption, we may find objects
$P,Q\in\Ob(\cC)$ with isomorphisms
$$
\omega_P:h_P\isom L_1:=\lim_{\Ob(I)}F\circ\mathrm{ob}
\qquad
\omega_Q:h_Q\isom L_2:=\lim_{\rMorph(I)}F\circ\mathrm{ob}\circ\mathrm{t}
\qquad
\text{in $\cC^\wedge_\sV$}
$$
for any universe $\sV$ such that $\cC$ has $\sV$-small
$\Hom$-sets. We define two morphisms of presheaves
$a,b:L_1\to L_2$ by the rules :
$$
a(\mu):=\mu*\mathrm{t}
\qquad\text{and}\qquad
b(\mu):=(F*\mathrm{m})\odot(\mu*\mathrm{s})
$$
for every $X\in\Ob(\cC)$ and every cone
$\mu:c_X\Rightarrow F\circ\mathrm{ob}$. Then $a$ and $b$
are represented by unique morphisms $\alpha,\beta:P\to Q$
in $\cC$ (see example \ref{rem_represent-morph}(i)) such that
$$
a\circ\omega_P=\omega_Q\circ h_\alpha
\qquad\text{and}\qquad
b\circ\omega_P=\omega_Q\circ h_\beta.
$$
Lastly, let $E:=\Equal(\alpha,\beta)$ (notation of example
\ref{ex_equalizers}(ii)), and notice that the tautological
cone for $E$ yields a monomorphism $E\to h_P$, so we may
regard $E$ as a subobject of $h_P$; it follows easily that
$\omega_P$ restricts to an isomorphism between $E\subset h_P$
and the presheaf $E'\subset L_1$ such that, for every
$X\in\Ob(\cC)$, the set $E'(X)$ consists of all the cones
$\mu:c_X\Rightarrow F\circ\mathrm{ob}$ with $a(\mu)=b(\mu)$.
A simple inspection shows that the latter condition describes
precisely the cones $c_X\Rightarrow F$; on the other hand, by
assumption $E$ is representable in $\cC$, so the same holds
for the limit of $F$.
\end{proof}

\begin{example}\label{ex_complete-cats}
(i)\ \
The category $\Set$ is complete and cocomplete. More precisely,
for any small family of sets $S_\bullet:=(S_i~|~i\in I)$ we have
natural representatives for the product and coproduct of
$S_\bullet$ : namely, for the product one may take the usual
cartesian product of the sets $S_i$, and for the coproduct
one may take the disjoint union
$$
\coprod_{i\in I}S_i:=\bigcup_{i\in I}\ (\{i\}\times S_i).
$$
Likewise, we have natural representatives for the equalizer
and coequalizer of any pair of maps $f,g:S\to S'$ : namely,
the equalizer of $f$ and $g$ is represented by the subset
$$
\Equal(f,g):=\{s\in S~|~f(s)=g(s)\}
$$
and the coequalizer is represented by the quotient $S'\!/\!\sim$,
where $\sim$ is the smallest equivalence relation on $S'$ whose
graph contains the subset $\{(f(s),g(s))~|~s\in S\}$.
Taking into account proposition \ref{prop_complete-criteria}(i,ii),
we get natural representatives for all limits and colimits in
$\Set$. For instance, by inspecting the proof of proposition
\ref{prop_complete-criteria}(i), we see that the colimit
of a functor $F:I\to\Set$ (for any small category $I$) is
represented by the quotient
$$
C(F):=\Bigl(\coprod_{i\in\Ob(I)}Fi\Bigr)/\!\sim
$$
where $\sim$ is the smallest equivalence relation such that
$(i,s)\sim (j,F(\phi)(s))$ for every morphism $\phi:i\to j$
in $I$ and every $s\in Fi$. In the special case where $I$
is filtered, this equivalence relation can be described
more explicitly : namely $(i,s)\sim(i',s')$ if and only if
there exist $j\in\Ob(I)$ and morphisms $\phi:i\to j$,
$\phi':i'\to j$ such that $F(\phi)(s)=F(\phi')(s')$.

For a general small category $I$, to every $i\in\Ob(I)$ we
may attach a natural map $\tau_i:Fi\to C(F)$ : namely,
$\tau_i(s)$ is the class of $(i,s)$ in $C(F)$, for every
$s\in Fi$; then the system $(\tau_i~|~i\in\Ob(I))$ yields
a universal cocone $\tau:F\Rightarrow c_{C(F)}$, in the
sense of remark \ref{rem_wishful}(iii).

(ii)\ \
In view of (i) and example \ref{ex_presheaves-on-sets}, we
can also give a more compact description of the limit of
any functor $F:I\to\Set$ : namely, if $\bone$ is any fixed
set with one element, the set $L(F):=\lim_IF(\bone)$,
consisting of all cones $c_\bone\Rightarrow F$, is another
natural representative for the limit of $F$. More precisely,
we get a natural isomorphism
$$
h_{L(F)}\isom\lim_IF
\qquad\text{in $\Set^\wedge$}
$$
as follows. A map of sets $S\to L(F)$ amounts to a system
of cones $\{\tau_s:c_\bone\Rightarrow F~|~s\in S\}$, and the
latter is naturally identified with a cone $c_S\Rightarrow F$
(details left to the reader). We easily deduce the following
natural universal cone $\tau:c_{L(F)}\Rightarrow F$. For every
$i\in\Ob(I)$, the map $\tau_i:L(F)\to Fi$ assigns to any
cone $\mu:c_\bone\Rightarrow F$ the image of $\mu_i:\bone\to Fi$.
\end{example}

\begin{example}\label{ex_big-complete}
(i)\ \
More generally, for every category $\cC$, the category
$\cC^\wedge$ is complete and cocomplete. Indeed, let $I$
be any small category, and $F:I\to\cC^\wedge$ any functor.
For every $X\in\Ob(\cC)$ we obtain a functor $F_X:I\to\Set$
such that $i\mapsto Fi(X)$ for every $i\in\Ob(I)$ and
$(\phi:i\to j)\mapsto(F\phi(X):Fi(X)\to Fj(X))$ for every
morphism $\phi$ of $I$; we let $LX$ be the natural
representative for $\lim_IF_X$, and $\tau^X:c_{LX}\Rightarrow F_X$
the corresponding universal cone, as described in example
\ref{ex_complete-cats}(ii). Every morphism $f:X\to Y$ in $\cC$
induces a natural transformation $F_f:F_Y\Rightarrow F_X$ by
the rule : $i\mapsto Fi(f):F_Y(i)\to F_X(i)$, whence a natural
transformation $\lim_IF_f:\lim_IF_Y\Rightarrow\lim_IF_X$ (remark
\ref{rem_wishful}(ii)) which is in turn represented by a
unique map of sets $Lf:LY\to LX$ such that
$$
\tau^X\odot c_{Lf}=F_f\odot\tau^Y
$$
(remark \ref{rem_represent-morph}(ii)). It is then easily seen
that the rules : $X\mapsto LX$ and $f\mapsto Lf$ yield a well
defined presheaf $L$ on $\cC$. Also, for every $i\in\Ob(I)$
the rule $X\mapsto\tau^X_i$ yields a morphism of presheaves
$\tau^\bullet_i:L\to Fi$, and the rule $i\mapsto\tau^\bullet_i$
yields a cone $\tau:c_L\Rightarrow F$.
Now, let $\eta:c_G\Rightarrow F$ be any cone with vertex
$G\in\Ob(\cC^\wedge)$. For every $X\in\Ob(\cC)$ we deduce a
cone $\eta^X:c_{G,X}\Rightarrow F_X$ by the rule :
$i\mapsto\eta_i(X):GX=G_X(i)\to F_X(i)$, whence a unique
morphism $t_X:GX\to LX$ such that $\tau^X\odot c_{t_X}=\eta^X$.
With this notation, a direct inspection shows that the
diagram
$$
\xymatrix{ c_{G,Y} \ar@{=>}[r]^-{\eta^Y} \ar@{=>}[d]_{c_{G,f}} &
F_Y \ar@{=>}[d]^{F_f} \\
c_{G,X} \ar@{=>}[r]^-{\eta^X} & F_X
}$$
commutes for every morphism $f:X\to Y$ of $\cC$. Thus,
$\tau^X\odot c_{t_X\circ G,f}=F_f\odot\tau^Y\odot c_{t_Y}=
\tau^X\odot c_{Lf\circ t_Y}$, whence $t_X\circ Gf=Lf\circ t_Y$,
by the universality of $\tau^X$. In other words, the rule :
$X\mapsto t_X$ yields a morphism $t:G\to L$ in $\cC^\wedge$,
and by construction, it is clear that $t$ is the unique such
morphism with  $\tau\odot c_t=\eta$. This shows that $\tau$
is a universal cone, whence the completeness of $\cC^\wedge$.
A similar argument proves the cocompleteness of $\cC^\wedge$ :
the details are left to the reader.

(ii)\ \
Furthermore, let $F:I\to\cC$ be a functor from a small
category $I$, whose limit $\cL$ (which is an object of
$\cC^\wedge_\sV$ for a suitable universe $\sV$) is representable
by some $L\in\Ob(\cC)$. Let $\tau:c_L\Rightarrow F$ be a
universal cone and denote by $\omega:h_L\isom\cL$ the
isomorphism corresponding to $\tau$, as in remark
\ref{rem_wishful}(iii).. Then the limit of
$h_\cC\circ F:I\to\cC^\wedge_\sV$ is represented by $h_L$ and
$h_\cC*\tau:c_{h_L}\Rightarrow h_\cC\circ F$ is a universal
cone. Indeed, let $t:c_\cL\Rightarrow h_\cC\circ F$
be the tautological cone (remark \ref{rem_wishful}(iv)); a
simple inspection shows that the resulting diagram
$$
\xymatrix{ & c_{h_L} \ar@{=>}[ld]_{c_\omega} \ar@{=>}[rd]^{h_\cC*\tau} \\
c_\cL \ar@{=>}[rr]^-{t} & & h_\cC\circ F
}$$
commutes. On the other hand, we know that $t$ is universal,
so the same follows for $h_\cC*\tau$.
\end{example}

\begin{example}\label{ex_cat-cats}
(i)\ \ 
Also the category $\bCat$ is complete and cocomplete. Indeed :

$\bullet$\ \
If $F:I\to\bCat$ is any functor from a small category $I$,
the limit of $F$ is representable by the category $\cC$
with $\Ob(\cC):=\lim_I\Ob\circ F$, where $\Ob:\bCat\to\Set$
is the functor of example \ref{ex_universe}(ii). Hence,
an object of $\cC$ is a family $(C_i~|~i\in\Ob(I))$ where
$C_i\in\Ob(\cC_i)$ for every $i\in I$, and $F(\phi)(C_i)=C_j$
for every $\phi:i\to j$ in $I$. For any two objects
$C_\bullet:=(C_i~|~i\in I)$ and $C'_\bullet:=(C'_i~|~i\in I)$,
notice that the rule $i\mapsto\Hom_{\cC_i}(C_i,C'_i)$ defines
a functor $H_{C_\bullet,C'_\bullet}:I\to\Set$; namely, to every
$\phi:i\to j$ in $I$ we assign the map
$H_{C_\bullet,C'_\bullet}(i)\to H_{C_\bullet,C'_\bullet}(j)$ given by
$F(\phi):\rMorph(\cC_i)\to\rMorph(\cC_j)$. Then we set
$$
\Hom_\cC(C_\bullet,C'_\bullet):=\lim_IH_{C_\bullet,C'_\bullet}.
$$
The composition of morphisms in $\cC$ is induced by the
composition laws of the categories $\cC_i$, in the obvious
way. The obvious projection functors $\cC\to\cC_i$ yield
a universal cone for $\lim_IF$.

$\bullet$\ \
Especially, for any pair of functors
$$
\cA\xrightarrow{\ F\ }\cC\xleftarrow{\ G\ }\cB
$$
the fibre product (in the category $\bCat$) of $F$ and $G$
is represented by the category :
$$
\cA\times_{(F,G)}\cB
$$
whose set of objects is $\Ob(\cA)\times_{\Ob(\cC)}\Ob(\cB)$; the
morphisms $(A,B)\to(A',B')$ are the pairs $(f,g)$ where $f:A\to A'$
(resp. $g:B\to B'$) is a morphism in $\cA$ (resp. in $\cB$), such
that $Ff=Gg$. If the notation is not ambiguous, we may also denote
this category by $\cA\times_\cC\cB$. In case $\cB$ is a subcategory
of $\cC$ and $G$ is the natural inclusion functor, we also write
$F^{-1}\cB$ instead of $\cA\times_\cC\cB$. The two obvious functors
$\cA\leftarrow\cA\times_\cC\cB\to\cB$ provide a universal cone.

$\bullet$\ \
If $(\cC_i~|~i\in I)$ is any small family of categories, the
coproduct $\cC:=\amalg_{i\in I}\cC_i$ is the category whose set
of objects is the disjoint union $\amalg_{i\in I}\Ob(\cC_i)$
(see example \ref{ex_complete-cats}(i)). Hence, an object of
$\cC$ is a pair $(i,C)$, with $i\in I$ and $C\in\Ob(\cC_i)$.
For any two such objects $(i,C)$ and $(i',C')$, we have
$$
\Hom_\cC((i,C),(i',C'))=\left\{\begin{array}{ll}
                        \Hom_{\cC_i}(C,C') & \text{if $i=i'$} \\
                        \emptyset & \text{otherwise}.
                              \end{array}\right.
$$
The obvious inclusion functors $\cC_i\to\cC$ provide a
universal cocone. The proof of the cocompleteness of
$\bCat$ shall be postponed to example \ref{ex_Cat-is-coco}
(see \cite[Prop.5.1.7]{Bor} for a more constructive proof,
and also example \ref{ex_fil-colim-in-Cat}).

(ii)\ \
For any category $\cA$, the objects of the fibre product
category $\sMorph(\cA)\times_{(\st,\ss)}\sMorph(\cA)$ are just
the composable pairs of morphisms of $\cA$ (notation of
\eqref{subsec_Morph-cat}), and the morphisms are the
pairs of commutative square diagrams that have one vertical
edge in common. Thus, we have a natural {\em composition functor}
$$
\sc:\sMorph(\cA)\times_{(\st,\ss)}\sMorph(\cA)\to\sMorph(\cA)
$$
such that $\sc(f,g):=g\circ f$ for every composable
pair of morphisms of $\cA$, and where $\sc$ assigns
to every pair of commutative square diagrams
$$
\xymatrix{
A \ar[r]^-f \ar[d]_{h_1} & B \ar[r]^-g \ar[d]_{h_2} &
C \ar[d]^{h_3} \\
A' \ar[r]^-{f'} & B' \ar[r]^-{g'} & C'
}$$
their composition, which is the square diagram
whose top (resp. bottom) horizontal arrow is $g\circ f$
(resp. $g'\circ f'$) and whose vertical arrows are
$h_1$ and $h_3$.

(iii)\ \
We can use the constructions of (ii) to rephrase the
definition of the composition laws for natural
transformations. Indeed, let us return to the situation
of \eqref{subsec_Godem-prod}, and denote by
$$
\tilde\alpha,\tilde\beta:\cA\to\rMorph(\cB)
$$
the functors attached to $\alpha$ and $\beta$ as in
\eqref{subsec_Morph-cat}. Clearly
$\st\circ\tilde\alpha=\ss\circ\tilde\beta$, whence
a unique functor
$$
(\tilde\alpha,\tilde\beta):\cA\to
\sMorph(\cB)\times_{(\st,\ss)}\sMorph(\cB)
$$
whose composition with the first (second) projection
to $\sMorph(\cB)$ equals $\tilde\alpha$ (resp.
$\tilde\beta$). Then it is easily seen that
$\sc\circ(\tilde\alpha,\tilde\beta):\cA\to\sMorph(\cA)$
is the (unique) functor associated to $\beta\odot\alpha$.

The Godement product of two natural transformations
$\alpha$ and $\alpha'$ as in \eqref{subsec_Godem-prod}
can be similarly described : to this aim, notice that
$$
(G'*\alpha)\odot(\alpha'*F)=\alpha'*\alpha=
(\alpha'*G)\odot(F'*\alpha)
$$
so we are reduced to exhibit the functors
$\cA\to\sMorph(\cC)$ associated to $G'*\alpha$ and
$\alpha'*F$. However, it is easily seen that the latter
are respectively $\sMorph(G')\circ\tilde\alpha$ and
$\tilde\alpha{}'\circ F$ (details left to the reader).
These considerations shall be amplified in section
\ref{sec_2Cats}, in order to construct analogous
composition laws for pseudo-natural transformations :
see remark \ref{rem_pseudo-natural}.
\end{example}

\begin{example}\label{ex_pre-misc}
Let $\cC$ be any category, and $\phi:X\to Y$ any morphism
of $\cC$. The {\em image} of $\phi$
$$
\Img(\phi)
$$
is defined the smallest of the subobjects $Y'\to Y$ such that
$\phi$ factors through a (necessarily unique) morphism $X\to Y'$.
If the image of $\phi$ exists, clearly a subobject $Z'\to Y$
contains $\Img(\phi)$ if and only if $\phi$ factors through $Z'$.

(i)\ \
In this generality, the image of $\phi$ does not necessarily
exists; however, suppose that $\cC$ is finitely complete
and the limit ({\em i.e.} the intersection) of any system
of subobjects of $Y$ is representable in $\cC$ (these limits
are well defined $\sU$-presheaves, for a sufficiently large
universe $\sU$). Then we claim that $\phi$ admits an image
$i:Z\to Y$. Indeed, let $\cF$ be the family of sub-objects
$Y'\to Y$ such that $\phi$ factors through $Y'$; choose
-- for every equivalence class $c\in\cF$ -- a representing
monomorphism $Y_c\to Y$, and let $I$ be the full subcategory
of $\cC\!/Y$ such that $\Ob(I)=\{Y_c\to Y~|~c\in\cF\}$. Then
$I$ is a $\sU$-small category, and we let $i$ (any
representing object for) the limit of the inclusion functor
$\iota:I\to\cC\!/Y$.
The resulting system of morphisms $(X\to Y_c~|~c\in\cF)$
is a cone on the indexing category $I$, whence a unique
morphism $p:X\to Z$, as required. To show that $p$ is
an epimorphism, consider any two morphisms $f,g:Z\to Z'$
(for any $Z'\in\Ob(\cC)$) such that $f\circ p=g\circ p$,
and let $E\in\Ob(\cC)$ be any object representing the
equalizer of $f$ and $g$. Thus, the induced morphism
$\psi:E\to Z$ a monomorphism, hence $E$ is a subobject of
of $Y$, and $\phi$ still factors through $E$. By the
minimality of $Z$, we deduce that $\psi$ is an isomorphism,
and therefore $f=g$, whence the claim.

(ii)\ \
Dually, if $\cC$ is finitely cocomplete and the colimit
of any system of quotients of $X$ is representable in
$\cC$, then the family of all quotients of $X$ through
which $\phi$ factors admits a minimal element $\psi':X\to Z$,
such that a quotient $Z'$ of $X$ maps to $Z$ if and only
if $\phi$ factors through $Z'$. Arguing as in (i), we see
that the resulting morphism $Z'\to Y$ is also a monomorphism.
We call $Z$ the {\em coimage} of $\phi$, and we denote it
$\mathrm{Coim}(\phi)$.

(iii)\ \
If both the image and the coimage of $\phi$ exist in $\cC$,
we get a natural factorization of $\phi$
$$
X\to\Img(\phi)\xrightarrow{\ \omega\ }\mathrm{Coim}(\phi)\to Y.
$$
But $\omega$ is not necessarily an isomorphism,
in this generality.
\end{example}

\begin{example}\label{ex_comma-adjunction}
(i)\ \
In the situation of \eqref{subsec_comma-adjunction}, let
$\eta_\bullet:\one_\cB\Rightarrow FG$ and
$\eps_\bullet:GF\Rightarrow\one_\cA$ be the unit and counit of
the adjunction $\theta$, and suppose that {\em all the fibre
products of $\cB$ are representable}; then for every
$B\in\Ob(\cB)$ the functor $G_{|B}:\cB/B\to\cA/GB$ admits a
right adjoint, denoted
$$
F_{|B}:\cA/GB\to\cB/B.
$$
Indeed, for every object $(g:A\to GB)$ of $\cA/GB$ let us fix
a cartesian diagram :
$$
\xymatrix{ F^*A \ar[r]^-{F^*g} \ar[d]_{\eta^*_A} & B \ar[d]^{\eta_B} \\
FA \ar[r]^-{Fg} & FGB.
}$$
If $h/GB:(g:A\to GB)\to(g':A'\to GB)$
is a morphism in $\cA/GB$, the universal property of the fibre
product yields a unique morphism $F^*h/B:F^*g\to F^*g'$ such
that $F^*g'\circ F^*h=F^*g$ and
$F^*g'\circ\eta^*_{A'}=\eta^*_A\circ Fh$. With this notation, we set
$$
F_{|B}(g):=F^*g
\qquad\text{and}\qquad
F_{|B}(h/GB):=F^*h/B
$$
for every such $g$ and $h/GB$. It is easily seen that these
rules define a functor as sought. In order to check that
$F_{|B}$ is right adjoint to $G_{|B}$, consider any morphism
$$
h/GB:G_{|B}(B'\xrightarrow{f}B)\to(A\xrightarrow{g}GB)
\qquad
\text{in $\cA/GB$}.
$$
Hence, $h:GB'\to A$ is a morphism in $\cA$ with $g\circ h=Gf$,
and we notice that :
$$
Fg\circ\theta(h)=Fg\circ Fh\circ\eta_{B'}=FGf\circ\eta_{B'}=
\eta_B\circ f.
$$
It follows that the pair $(\theta_{AB'}(h),f)$ determines a unique
morphism $k:B'\to F^*A$ such that $\eta^*_A\circ k=\theta_{AB'}(h)$
and $F^*g\circ k=f$. With this notation, we set
$(\theta_{|B})_{g,f}(h/GB):=k/B:f\to F^*g$ in $\cB/B$. Conversely,
to every morphism $k/B:(f:B'\to B)\to F_{|B}(g:A\to GB)$ in $\cB/B$,
we attach the morphism $h:=\theta^{-1}_{AB'}(\eta^*_A\circ k):GB'\to A$.
Let us show that $h/GB:Gf\to g$ is a morphism in $\cA/GB$;
recalling that $F^*g\circ k=f$, it suffices to compute :
$$
g\circ h=g\circ\eps_A\circ G(\eta^*_A\circ k)=
\eps_{GB}\circ GFg\circ G(\eta^*_A\circ k)=
\eps_{GB}\circ G(\eta_B\circ F^*g\circ k)=Gf
$$
where the last equality follows from the triangular identities
of \eqref{subsec_adj-pair}. Hence the map
$$
(\theta_{|B})_{g,f}:\Hom_{\cA/GB}(G_{|B}f,g)\to\Hom_{\cB/B}(f,F_{|B}g)
$$
is a bijection with inverse given by the rule :
$k/B\mapsto\theta^{-1}_{AB'}(\eta^*_A\circ k)/GB$. Indeed,
by definition
$(\theta_{|B})_{g,f}(\theta^{-1}_{AB'}(\eta^*_A\circ k)/GB)$ is
the morphism $f\to F_{|B}g$ of $\cB/B$ determined by the pair
$(\theta_{AB'}\theta^{-1}_{AB'}(\eta^*_A\circ k),f)=
(\eta^*_A\circ k),f)$, which is just $k/B$, and on the other
hand, $\theta^{-1}_{AB'}(\eta^*_A\circ(\theta_{|B})_{g,f}(h/GB))=
\theta^{-1}_{AB'}\theta_{AB'}(h)/GB=h/GB$, whence the contention.
The naturality of $(\theta_{|B})_{g,f}$ with respect to $g$
and $f$ follows by a simple inspection : details left to the
reader.

(ii)\ \
Dually, if {\em all the coproducts of $\cA$ are representable},
then the functor ${}_{A|}F:A/\cA\to FA/\cB$ admits a left
adjoint ${}_{A|}G:FA/\cB\to A/\cA$, which the reader is
invited to spell out.

(iii)\ \
Keep the situation of (i), and let $F':\cB\to\cC$ be
another functor with a left adjoint $G':\cC\to\cB$.
Then we have an isomorphism of functors :
$$
F'_{|C}\circ F_{|G'C}\isom(F'F)_{|C}
\qquad
\text{for every $C\in\Ob(\cC)$}.
$$
Indeed, clearly we have $G_{|G'C}\circ G'_C=(GG')_{|C}$, so the
assertion follows from example \ref{ex_presheaves-on-sets}(ii).
\end{example}

\subsection{Adjunctions and Kan extensions}\label{sec_Fubini}
Let $\phi:I\to J$ be a functor between small categories.
The calculation of limits indexed by $I$ can be sometimes
simplified by a Fubini-style technique of ``integration
along the fibres of $\phi$'', which replaces any given
functor $F:I\to\cC$ (to any category $\cC$) with another
functor $J\to\cC^\wedge$ whose limit is isomorphic to that
of $F$. Namely, let $\sV$ be a universe containing $\sU$,
such that $\cC$ has $\sV$-small $\Hom$-sets; notice first
that by virtue of remark \ref{rem_wishful}(ii), the limit
construction yields a well defined functor
$$
\lim_I:\bFun(I,\cC)\to\cC^\wedge_\sV
\qquad
F\mapsto\lim_IF.
$$
Moreover, if $\cC$ is complete (resp. if $I$ is finite
and $\cC$ is finitely complete), this functor is naturally
isomorphic to the composition of the Yoneda embedding
$h_\cC$ and a functor
$$
\Lim_I:\bFun(I,\cC)\to\cC
$$
which is right adjoint to the {\em constant functor}
\set\begin{equation}\label{eq_const-functor}
c:\cC\to\bFun(I,\cC)
\qquad
X\mapsto c_X
\qquad
(f:X\to Y)\mapsto(c_f:c_X\Rightarrow c_Y)
\end{equation}
(notation of \eqref{subsec_wishful}), and especially, it is
well defined up to isomorphism (example
\ref{ex_presheaves-on-sets}(ii)). We define now as follows
a functor
$$
\int^\wedge_\phi:\bFun(I,\cC)\to\bFun(J,\cC^\wedge_\sV).
$$
For every functor $F:I\to\cC$, every $j\in\Ob(J)$,
and every morphism $g:j\to j'$ in $J$ we set
$$
\int^\wedge_\phi F(j):=\lim_{j/\phi I}F\circ\st_j
\qquad
\int^\wedge_\phi F(g):=\lim_{g/\phi I}\one_\cC
$$
(notation of \eqref{subsec_fibreovercat} and remark
\ref{rem_wishful}(i)). We may then state :

\begin{proposition}\label{prop_Fubini}
With the notation of \eqref{sec_Fubini}, the diagram
$$
\xymatrix{
\bFun(I,\cC) \ar[rr]^-{\lim_I} \ar[d]_{\int^\wedge_\phi}
& & \cC^\wedge_\sV \\
\bFun(J,\cC^\wedge_\sV) \ar[rru]_{\Lim_J}
}$$
is {\em essentially commutative}, {\em i.e.} there is a
natural isomorphism of functors :
$\lim_I\isom\Lim_J\circ\int^\wedge_\phi$.
\end{proposition}
\begin{proof} Notice that the functor $\Lim_J$ is well defined
by virtue of example \ref{ex_big-complete}. Let $F:I\to\cC$ be
any functor; unwinding the definitions, we see that
$\Lim_J\circ\int^\wedge_\phi F$ is isomorphic to the presheaf
that assigns to every $X\in\Ob(\cC)$ the set of all
compatible systems
$$
\tau^\bullet_\bullet:=
(\tau^j_\bullet:c_X\Rightarrow F\circ\st_j~|~j\in\Ob(J))
$$
such that $\tau^{j'}_\bullet=\tau^j_\bullet*(g/\phi I)$
for every morphism $g:j\to j'$ in $J$. Explicitly,
$\tau^\bullet_\bullet$ is the datum, for every
$(h:j'\to\phi i)\in\Ob(J/\phi I)$, of a morphism
$\tau^j_h:X\to Fi$ in $\cC$, such that
$$
\tau^j_{h\circ g}=\tau^{j'}_h
\quad\text{and}\quad
Ff\circ\tau^j_h=\tau^j_{\phi(f)\circ h}
$$
for every $g$ as above, and every morphism $f:i\to i'$
in $I$. Clearly, every such system is determined by the
family $(\bar\tau_i:=\tau^{\phi i}_{\one_{\phi i}}~|~i\in\Ob(I))$,
and the latter is a cone $\bar\tau_\bullet:c_X\Rightarrow F$.
Conversely, given any cone $\mu_\bullet:c_X\Rightarrow F$, the
family $\bar\mu{}^\bullet_\bullet:=(\mu_\bullet*\st_j~|~j\in\Ob(J))$
is a compatible system of the foregoing type. A simple
inspection shows that the rules
$\tau_\bullet^\bullet\mapsto\bar\tau_\bullet$ and
$\mu_\bullet\mapsto\bar\mu{}^\bullet_\bullet$ are mutually
inverse, whence the proposition.
\end{proof}

\begin{remark}\label{rem_Lims-and-Colims}
(i)\ \
In the situation of \eqref{sec_Fubini}, suppose
moreover that $\cC$ is complete. Then $\int^\wedge_\phi$
is isomorphic to the composition of the functor $\bFun(J,h_\cC)$
(notation of remark \ref{rem_opposite-Fun}(ii) and
\eqref{subsec_yoneda}) and a functor
$$
\int_\phi:\bFun(I,\cC)\to\bFun(J,\cC)
$$
called the {\em right Kan extension along $\phi$}. Taking into
account example \ref{ex_big-complete}(ii), proposition
\ref{prop_Fubini} implies that the resulting diagram
$$
\xymatrix{
\bFun(I,\cC) \ar[rr]^-{\Lim_I} \ar[d]_{\int_\phi} & & \cC \\
\bFun(J,\cC) \ar[rru]_{\Lim_J}
}$$
is also essentially commutative (details left to the reader).

(ii)\ \
By remark \ref{rem_wishful}(ii), we may also define a colimit
functor
$$
\colim_I:\bFun(I,\cC)\to\cC^{o\wedge o}_\sV
\qquad
F\mapsto\colim_IF
$$
and notice that the isomorphism of remark
\ref{rem_opposite-Fun}(i) identifies this functor with
the opposite of the functor $\lim_{I^o}$ defined on
$\bFun(I^o,\cC^o_\sV)$. Consequently, if $\cC$ is cocomplete
(resp. if $I$ is finite and $\cC$ is finitely cocomplete),
$\colim_I$ is isomorphic to the composition of $h^o_{\cC^o}$
and of a functor
$$
\Colim_I:\bFun(I,\cC)\to\cC
$$
which is likewise identified with the opposite of
the functor $\Lim_{I^o}:\bFun(I^o,\cC^o)\to\cC^o$. Notice
as well that $\Colim_I$ is left adjoint to the constant
functor \eqref{eq_const-functor}.
Moreover, the opposite of the functor $\int^\wedge_{\phi^o}$
on $\bFun(I^o,\cC^o)$ is identified with a functor
$$
\int^\phi_\wedge:\bFun(I,\cC)\to\bFun(J,\cC^{o\wedge o}_\sV).
$$
Explicitly, we have
$$
\int_\wedge^\phi F(j)=\colim_{\phi I/j}F\circ\ss_j
\qquad
\text{for every functor $F:I\to\cC$ and every $j\in\Ob(J)$}
$$
(where $\ss_j$ is the source functor) and we get an
essentially commutative diagram :
$$
\xymatrix{
\bFun(I,\cC) \ar[rr]^-{\colim_I} \ar[d]_{\int_\wedge^\phi}
& & \cC^{o\wedge o}_\sV \\
\bFun(J,\cC^{o\wedge o}_\sV) \ar[rru]_{\Colim_J}
}$$
Notice the natural identification :
$\cC^{o\wedge o}_\sV\isom\bFun(\cC^o,\sV\tdu\Set^o)$, again due
to remark \ref{rem_opposite-Fun}(i). If $\cC$ is cocomplete,
the opposite of the functor $\int_{\phi^o}$ is naturally
identified with a functor
$$
\int^\phi:\bFun(I,\cC)\to\bFun(J,\cC)
$$
called the {\em left Kan extension along $\phi$}, fitting
into the essentially commutative diagram :
$$
\xymatrix{
\bFun(I,\cC) \ar[rr]^-{\Colim_I} \ar[d]_{\int^\phi}
& & \cC \\
\bFun(J,\cC). \ar[rru]_{\Colim_J}
}$$

(iii)\ \
Notice that the functors $\int^\wedge_\phi$ and $\int^\phi_\wedge$
are defined more generally whenever $I$ is small and $J$ has
small $\Hom$-sets. Likewise, if $\cC$ is complete (resp.
cocomplete) the right (resp. left) Kan extension along $\phi$
is well defined under the same weakened assumptions.

(iv)\ \
In the same vein, if the category  $j/\phi I$ (resp. $\phi I/j$)
is finite for every $j\in\Ob(J)$, for the right (resp. left) Kan
extension along $\phi$ to be well defined it suffices that $\cC$
is finitely complete (resp. finitely cocomplete).
\end{remark}

\begin{theorem}\label{th_Kan-ext}
Let $I$, $J$ and $\cC$ be categories fulfilling the conditions
of either part {\em (iii)} or part {\em (iv)} of remark
{\em\ref{rem_Lims-and-Colims}}, and $\phi:I\to J$ any functor.
Then the right (resp. left) Kan extension along $\phi:I\to J$
is right (resp. left) adjoint to the functor
$$
\bFun(\phi,\cC):\bFun(J,\cC)\to\bFun(I,\cC)
\qquad
G\mapsto G\circ\phi.
$$
\end{theorem}
\begin{proof} (Cp. \cite[Th.3.7.2]{Bor} or \cite[Th.2.3.3]{Ka-Sch}).
Pick a universe $\sV$ containing $\sU$, and such that $\cC$
has $\sV$-small $\Hom$-sets. In light of Yoneda's lemma (proposition
\ref{prop_yoneda}(ii)) we get for every pair of functors $F:I\to\cC$
and $G:J\to\cC$ a natural bijection
$$
\Hom_{\bFun(J,\cC)}\Bigl(G,\int_\phi F\Bigr)\isom
\Hom_{\bFun(J,\cC^\wedge_\sV)}\Bigl(h_\cC\circ G,\int^\wedge_\phi F\Bigr)
$$
where $h_\cC:\cC\to\cC^\wedge_\sV$ is the Yoneda embedding of
\eqref{subsec_yoneda}. Hence, in order to prove the assertion for
right Kan extensions, it suffices to exhibit a natural bijection :
\set\begin{equation}\label{eq_Kan-extend}
\Hom_{\bFun(I,\cC)}(G\circ\phi,F)\isom
\Hom_{\bFun(J,\cC^\wedge_\sV)}\Bigl(h_\cC\circ G,\int^\wedge_\phi F\Bigr).
\end{equation}
Notice that, by virtue of proposition \ref{prop_yoneda}(ii), if
$X\in\Ob(\cC)$ and $j\in\Ob(J)$ are any two objects, a morphism
$h_X\to\int^\wedge_\phi F(j)$ in $\cC^\wedge_\sV$ is the same as the
datum of a cone $c_X\Rightarrow F\circ\st_j$, where
$\st_j:j/\phi I\to I$ is the target functor. Now, let
$\tau:G\circ\phi\Rightarrow F$ be any natural transformation;
for every $j\in\Ob(J)$ and every
$(i,f:j\to\phi(i))\in\Ob(j/\phi I)$ we set
$$
\tau'_{j,(i,f)}:=\tau_i\circ Gf:Gj\to Fi.
$$
It is easily seen that, for fixed $j\in\Ob(J)$, the system
$(\tau'_{j,(i,f)}~|~(i,f)\in\Ob(j/\phi I))$ gives a cone
$\tau'_j:c_{Gj}\Rightarrow F\circ\st_X$, {\em i.e.} a morphism
$h_{Gj}\to\int^\wedge_\phi F(j)$ in $\cC^\wedge_\sV$, and the rule
$j\mapsto\tau'_j$ for every $j\in\Ob(J)$ yields a natural
transformation $\tau':h_\cC\circ G\Rightarrow\int^\wedge_\phi F$.
Conversely, if $\mu:h_\cC\circ G\Rightarrow\int^\wedge_\phi F$ is
any such transformation, every $j\in\Ob(J)$ yields a cone
$\mu_j:c_{Gj}\Rightarrow F\circ\st_j$, and we define a
natural transformation $\mu':G\circ\phi\Rightarrow F$ by setting
$$
\mu'_i:=\mu_{\phi(i),(i,\one_{\phi(i)})}:G\circ\phi(i)\to Fi
\qquad
\text{for every $i\in\Ob(I)$}.
$$
The reader can check that $(\tau')'=\tau$ and $(\mu')'=\mu$
for every $\tau$ and $\mu$ as above, and clearly the resulting
bijections \eqref{eq_Kan-extend} are functorial in both $F$
and $G$. The assertion for left Kan extensions, in case $\cC$
is cocomplete, is dual to the foregoing, by virtue of remark
\ref{rem_opposite-Fun}(i).
\end{proof}

\begin{remark}\label{rem_was-cofinal}
(i)\ \
Let $F:\cB\to\cC$ be a functor from a small category $\cB$
to any category $\cC$ with small $\Hom$-sets; from theorem
\ref{th_Kan-ext} we obtain both a left and a right adjoint
for $F^\wedge_\sU=\bFun(F^o,\sU\tdu\Set)$, denoted respectively
$$
F_{\sU!}:\cB^\wedge_\sU\to\cC^\wedge_\sU
\qquad \text{and} \qquad
F_{\sU*}:\cB^\wedge_\sU\to\cC^\wedge_\sU.
$$
As usual, we drop the subscript $\sU$, unless the omission
may cause ambiguities.

(ii)\ \
With the notation of (i), notice that the diagram of
functors :
\set\begin{equation}\label{eq_ess-comm-dig}
{\diagram
\cB \ar[r]^-{h_\cB} \ar[d]_F & \cB^\wedge \ar[d]^{F_!} \\
\cC \ar[r]^-{h_\cC} & \cC^\wedge
\enddiagram}
\end{equation}
(whose horizontal arrows are the Yoneda embeddings) is
{\em essentially commutative}, {\em i.e.} the two compositions
$F_!\circ h_\cB$ and $h_\cC\circ F$ are isomorphic functors.
Indeed -- by proposition \ref{prop_yoneda}(ii) -- for every
$B\in\Ob(\cB)$, the objects $F_!h_B$ and $h_{FB}$ both represent
the functor
$$
\cC^\wedge\to\Set
\quad : \quad
\phi\mapsto\phi(FB)
$$
(the latter is a $\sV$-presheaf on $\cC^{\wedge o}$, for any
universe $\sV$ such that $\cC^\wedge$ has small $\Hom$-sets).

(iii)\ \
Moreover, taking into account example \ref{ex_complete-cats}
we see that there are natural choices for the left and right
adjoints of the functor $F^\wedge$. Namely :

$\bullet$\ \
For any presheaf $f$ on $\cB$, the presheaf $F_!f$ on $\cC$
can be given canonically by the rule :
\set\begin{equation}\label{eq_describe_f_lower}
C\mapsto\colim_{(C/F\cB)^o} f\circ\st^o_C
\qquad
\text{for every $C\in\Ob(\cC)$}
\end{equation}
where the colimits in this expression denote (by abuse of
notation) the natural representatives for the corresponding
presheaves, as described in example \ref{ex_complete-cats}(i).
Therefore, every element of $F_!f(C)$ is the equivalence class
of a pair $(\psi:C\to FB,s)$, where $\psi$ is any object of
$C/F\cB$ and $s\in fB$ is any element. If $\phi:C'\to C$ is
any morphism of $\cC$, the induced map
$F_!f(\phi):F_!f(C)\to F_!f(C')$ assigns to any such pair
the equivalence class of the pair $(\psi\circ\phi:C'\to FB,s)$.
Furthermore, the adjunction for the pair $(F_!,F^\wedge)$
provided by the proof of theorem \ref{th_Kan-ext} comes down
to the natural bijection
$$
\Hom_{\cB^\wedge}(f,F^\wedge g)\isom\Hom_{\cC^\wedge}(F_!f,g)
\qquad
\text{for every $f\in\Ob(\cB^\wedge)$ and $g\in\Ob(\cC^\wedge)$}
$$
which, to any morphism $\tau:f\to F^\wedge g$ in $\cB^\wedge$ and
to any $C\in\Ob(\cC)$ assigns the map of sets $F_!f(C)\to gC$
given by the rule :
$$
(\psi:C\to FB,s)\mapsto g(\psi)\circ\tau_B(s)
\qquad
\text{for every $(B,\psi)\in\Ob(C/F\cB)$ and every $s\in fB$}.
$$
The inverse to this bijection assigns to any morphism
$\mu:F_!f\to g$ and any $B\in\Ob(\cB)$ the map
$f(B)\to g(FB)$ given by the rule :
$$
s\mapsto\mu_{FB}(\one_{FB},s)
\qquad
\text{for every $s\in f(B)$}.
$$

$\bullet$\ \
And the presheaf $F_*f$ can be chosen canonically by the rule :
$$
C\mapsto\lim_{(F\cB/C)^o} f\circ\ss^o_C(\bone)
\qquad
\phi\mapsto\lim_{(F\cB/\phi)^o} \one_\cB(\bone)
$$
for every $C\in\Ob(\cC)$ and every morphism $\phi:C'\to C$ in
$\cC$ (notation of \eqref{subsec_fibreovercat}, and $\bone$ is
any fixed set with one element). More explicitly, under this
identification, the map $F_*f(\phi):F_*f(C)\to F_*f(C')$ is
given by the rule :
$$
(\tau:c_\bone\Rightarrow f\circ\ss^o_C)\mapsto
(\tau*(F\cB/\phi)^o:c_\bone\Rightarrow f\circ\ss^o_{C'}).
$$

(iv)\ \
Let $\sV$ be a universe such that $\sU\subset\sV$. Notice
that, with the canonical choices of (iii), we get a commutative
diagram of categories
$$
\xymatrix{
\cC^\wedge_\sU \ar[d] &
\cB^\wedge_\sU \ar[r]^-{F_{\sU!}} \ar[l]_-{F_{\sU*}} \ar[d] &
\cC^\wedge_\sU \ar[d] \\
\cC^\wedge_\sV & \cB^\wedge_\sV \ar[r]^-{F_{\sV!}} \ar[l]_-{F_{\sV*}}
& \cC^\wedge_\sV
}$$
whose vertical arrows are the inclusion functors.

(v)\ \
Lastly, with the canonical choices in (iii), we can make (ii)
more precise : we get a natural identification
$$
F_!h_B\isom h_{FB}
\qquad
\text{for every $B\in\Ob(\cB)$}.
$$
Namely, to any pair $(\psi':Y\to FB',s:B'\to B)$ consisting
of an object of $Y/F\cB$ and an element $s\in h_B(B')$, we
attach the element $\beta(\psi',s):=F(s)\circ\psi'\in h_{FB}(Y)$,
and it is easily seen that $\beta(\psi',s)$ depends only on
the class of $(\psi',s)$ in $F_!h_B(Y)$, whence the required
functorial bijection (details left to the reader). Notice
as well that -- under these identifications -- the unit of
the adjunction for the pair $(F_!,F^\wedge)$ explicited in
(iii) becomes the natural morphism
$$
h_B\to F^\wedge h_B
\quad :\quad
(s:B'\to B)\mapsto(Fs:FB'\to FB)
\qquad
\text{for every $B\in\Ob(B)$}.
$$
\end{remark}

\begin{example}\label{ex_was-lim-yoneda}
Let $\bone_\cB$ be a final object of $\cB^\wedge$ (see
example \ref{ex_equalizers}(v)); by specializing the
explicit description of $F_!$ provided by remark
\ref{rem_was-cofinal}(iii), we get a natural isomorphism
$$
F_!(\bone_\cB)\isom\pi_0(\bullet/F\cB)
\qquad\text{in $\cC^\wedge$}
$$
where $\pi_0(\bullet/F\cB)$ is the presheaf that assigns
to any $C\in\Ob(\cC)$ the set $\pi_0(C/F\cB)$, and to
any morphism $\phi:C'\to C$ in $\cC$ the map
$\pi_0(\phi/F\cB):\pi_0(C/F\cB)\to\pi_0(C'/F\cB)$
induced by $\phi/F\cB$ in the obvious way (notation of
\eqref{subsec_fibreovercat} : we leave the verification
to the reader).
\end{example}

\sset\subsubsection{}
For many questions, it is useful to know whether a given
functor commutes with a certain limit, or more generally
with a prescribed class of limits (or colimits), in the
sense explained by the following :

\begin{definition}\label{def_commute-with-lim}
Let $I$, $\cC$ and $\cD$ be any three categories, and
$F:I\to\cC$, $f:\cC\to\cD$ two functors. Pick also a
universe $\sV$ such that $I$, $\cC$ and $\cD$ are
$\sV$-small. 
\begin{enumerate}
\item
According to remark \ref{rem_wishful}(i) we have a
natural morphism in $\cC^\wedge_\sV$
\set\begin{equation}\label{eq_before-adjoint}
\lim_IF\to f^\wedge_\sV\lim_I\,(f\circ F)
\end{equation}
whence -- by the explicit adjunction of remark
\ref{rem_was-cofinal}(iii) -- a morphism in
$\cD^\wedge_\sV$
\set\begin{equation}\label{eq_F-commutes-with-f}
f_{\sV!}(\lim_IF)\to\lim_I\,(f\circ F).
\end{equation}
Then, we say that {\em $f$ commutes with the limit of $F$},
if \eqref{eq_F-commutes-with-f} is an isomorphism.
\item
We say that {\em $f$ commutes with the colimit of $F$},
if the dual condition holds, {\em i.e.} if $f^o$ commutes
with the limit of $F^o$, so the natural morphism
$$
f^o_{\sV!}(\colim_IF)\to\colim_I(f\circ F)
$$
is an isomorphism.
\item
We say that the functor $f$ is {\em left exact}, if $f$
commutes with all finite limits. Dually, we say that $f$
is {\em right exact\/} if it commutes with all finite
colimits. Finally, we say that $f$ is {\em exact\/} if
it is both left and right exact.
\end{enumerate}
It follows easily from remark \eqref{rem_was-cofinal}(iv)
that these definitions depend only on $f$ and $F$, and not
on the choice of the universe $\sV$.
\end{definition}

\begin{remark}\label{rem_representable-case}
(i)\ \
Keep the situation of definition \ref{def_commute-with-lim}.
The most interesting case is that where the limit of $F$
is representable, say by an object $L$ of $\cC$, so we
have an isomorphism $\beta:h_L\isom\lim_IF$. Indeed, in
this case, by remark \ref{rem_was-cofinal}(v) we deduce
a morphism
$$
\gamma:h_{fL}\isom f_!h_L\isom f_!\lim_IF\to\lim_I(f\circ F)
$$
that is an isomorphism if and only if $f$ commutes with
the limit of $F$, in which case we see that the limit of
$f\circ F$ is representable by $fL$. More precisely, let
$\tau:=\beta_L(\one_L):c_L\Rightarrow F$ be the universal
cone deduced from $\beta$; we claim that
$$
\gamma_{fL}(\one_{fL})=f*\tau:c_{fL}\Rightarrow f\circ F
$$
so $f$ commutes with the limit of $F$ if and only if
$f*\tau$ is a universal cone for $f\circ F$. Indeed,
notice that the composition
$h_L\to f^\wedge_\sV\lim_I(f\circ F)$ of $\beta$ and
\eqref{eq_before-adjoint} assigns, to any $C\in\Ob(\cC)$
and any morphism $s:C\to L$ in $\cC$, the cone
$(f*\tau)\odot c_{fs}:c_{fC}\Rightarrow f\circ F$. In light
of remark \ref{rem_was-cofinal}(iii), it follows that the
adjoint morphism $f_!h_L\to\lim_If\circ F$ assigns, to every
$D\in\Ob(\cD)$ and every pair $(\psi:D\to fC,s:C\to L)$ the
cone $(f*\tau)\odot c_{fs}\odot c_\psi:c_D\Rightarrow f\circ F$.
Lastly, notice that the isomorphism $f_!h_L\isom h_{fL}$
provided from remark \ref{rem_was-cofinal}(v) maps the pair
$(\one_{fL},\one_L)$ to $\one_{fL}\in h_{fL}(fL)$; we conclude
that $\gamma_{fL}(\one_{fL})=
(f*\tau)\odot c_{f\one_L}\odot c_{\one_{fL}}=f*\tau$, as claimed.

(ii)\ \
Dually, suppose that the colimit of $F$ is representable
by an object $C^o$ of $\cC^o$, and pick any universal
cocone $\mu:F\Rightarrow c_C$. Then $f$ commutes with
the colimit of $F$ if and only if $f*\mu$ is a universal
cocone for $f\circ F$.
\end{remark}

\begin{example}\label{ex_lims-and-representables}
(i)\ \
Let $\cC$ be any category, $\sV$ a universe such that $\cC$
has $\sV$-small $\Hom$-sets, and $f:\cC^o\to\sV\tdu\Set$ a
representable functor. Then it is easily seen that for every
small category $I$, every functor $F:I\to\cC$ and every universal
cocone $\tau:F\Rightarrow c_X$, the induced cone
$f*\tau^o:c_{fX}\Rightarrow f\circ F^o$ is universal. In view of
remark \ref{rem_representable-case}(i), we conclude that $f$
commutes with all the representable limits of $\cC^o$ (which are
the representable colimits of $\cC$). On the other hand, if
$\cC$ has small $\Hom$-sets and all small products of $\cC$
are representable, then $f$ commutes even with non-representable
colimits of $\cC$ : see example \ref{ex_rep-funct-commute-with-lims}(i).

(ii)\ \
Dually, if $f:\cC\to\Set$ is a functor representable by
some object $Y^o$ of the opposite category $\cC^o$, then
$f$ commutes with all representable limits of $\cC$, and
even with non-representable limits, if all small coproducts
of $\cC$ are representable and $\cC$ has small $\Hom$-sets.
\end{example}

\begin{example}\label{ex_cofiltered-comma}
(i)\ \
Let $\cA$ be any finitely complete category, and
$F:\cA\to\cB$ any left exact functor. Then the
category $B/F\cA$ is cofiltered, for every
$B\in\Ob(\cB)$. Indeed, if $\phi_i:B\to FA_i$ ($i=1,2$)
are any two objects of $B/F\cA$, the product $A_1\times A_2$
is representable by some object $A$ of $\cA$, and
a universal cone is given by a pair of morphisms
$(p_i:A\to A_i~|~i=1,2)$; by remark
\ref{rem_representable-case}(i), it follows that
the pair $(Fp_i:FA\to FA_i~|~i=1,2)$ is still a
universal cone for the product $FA_1\times FA_2$,
whence a unique morphism $\phi:B\to FA$ such that
$Fp_i\circ\phi=\phi_i$ for $i=1,2$. This shows that
$B/F\cA$ is codirected. Moreover, if $A$ is any final
object of $\cA$, then $FA$ is a final object of $\cB$,
hence $\Ob(B/F\cA)$ is non-empty. By remark \ref{rem_cofinal}(i),
it then remains only to check that $B/F\cA$ satisfies
the equalizing condition dual to that of definition
\ref{def_filtered-cols}(v). Namely, let $\phi_1$ and
$\phi_2$ as in the foregoing, and suppose we have a
pair of morphisms $\psi,\psi':A_1\to A_2$ such that
$F\psi\circ\phi_1=\phi_2=F\psi'\circ\phi_1$. The
equalizer of $\psi$ and $\psi'$ is representable
by some object $E$ of $\cA$, and a universal cone
for $E$ is given by a morphism $\beta:E\to A_1$
such that $\psi\circ\beta=\psi'\circ\beta$; by
remark \ref{rem_representable-case}(i), $FE$
represents the equalizer of $F\psi$ and $F\psi'$,
and $F\beta:FE\to FA_1$ still yields a universal cone.
There follows a unique morphism $\gamma:B\to FE$ in
$\cB$ such that $F\beta\circ\gamma=\phi_1$, whence
the claim.

(ii)\ \
Dually, if $\cA$ is finitely cocomplete and $F$ is
right exact, then the category $F\cA/B$ is filtered
for every $B\in\Ob(\cB)$.
\end{example}

\begin{proposition}\label{prop_double-double}
Let $I$, $\cC$, $\cD$ and $\cE$ be four categories, and
$F:I\to\cC$, $f:\cC\to\cD$ and $g:\cD\to\cE$ any three
functors. The following holds :
\begin{enumerate}
\item
If $f$ commutes with the limit of $F$, and $g$ commutes with
the limit of $f\circ F$, then $g\circ f$ commutes with the
limit of $F$.
\item
Especially, if both $f$ and $g$ are left (resp. right) exact,
then the same holds for $g\circ f$.
\end{enumerate}
\end{proposition}
\begin{proof} Clearly it suffices to show (i), and after replacing
$\sU$ by a larger universe, we may assume that the four categories
of the proposition are small. Set $L:=\lim_IF$, $L':=\lim_I(f\circ F)$,
$L'':=\lim_I(g\circ f\circ F)$, and denote
$$
\omega_{f,F}:L\to f^\wedge L'
\qquad\text{and}\qquad
\omega_{g,f\circ F}:L'\to g^\wedge L''
$$
the morphisms provided by remark \ref{rem_wishful}(i). Let
also $(\eps^f,\eta^f)$ (resp. $(\eps^g,\eta^g)$) be the units
and counits associated with the explicit adjunction $\theta_f$
(resp $\theta_g$) for the pair $(f_!,f^\wedge)$ (resp.
$(g_!,g^\wedge)$) provided by remark \ref{rem_was-cofinal}(iii);
by assumption, the adjoint morphisms
$$
\omega^*_{f,F}:=\eps^f_{L'}\circ(f_!\omega_{f,F}):f_!L\to L'
\qquad\text{and}\qquad
\omega^*_{g,f\circ F}:=\eps^g_{L''}\circ(g_!\omega_{g,f\circ F}):g_!L'\to L''
$$
are isomorphisms, so the same holds for $\omega^{**}_{g\circ f,F}:=
\omega^*_{g,f\circ F}\circ g_!\omega^*_{f,F}:g_!\circ f_!L\to L''$.
However :
$$
\omega^{**}_{g\circ f,F}=\eps^g_{L''}\circ
g_!(\omega_{g,f\circ F}\circ\omega^*_{f,F})
$$
corresponds, under the adjunction $\theta_g$, to the morphism
$\omega_{g,f\circ F}\circ\omega^*_{f,F}:f_!L\to g^\wedge L''$.
The latter, in turns, corresponds -- under the adjunction
$\theta_f$ -- to the morphism
$$
\omega_{g\circ f,F}:=
f^\wedge(\omega_{g,f\circ F}\circ\omega^*_{f,F})\circ\eta^f_L:
L\to f^\wedge\circ g^\wedge L''=(g\circ f)^\wedge L''
$$
and notice that $f^\wedge(\omega^*_{f,F})\circ\eta^f_L=\omega_{f,F}$.
Taking into account \eqref{eq_respect}, we conclude that
$\omega_{g\circ f,F}$ is precisely the natural morphism denoted
$\lim_{\one_I}(g\circ f)$ in remark \ref{rem_wishful}(i).

Summing up, we have shown that, under the adjunction
$\theta_g\circ\theta_f$ for the pair
$(g_!\circ f_!,(g\circ f)^\wedge)$ (notation of remark
\ref{rem_adjoint-transf}(i)), the morphism $\omega_{g\circ f,F}$
corresponds to the isomorphism $\omega^{**}_{g\circ f,F}$.
Then the assertion follows from example
\ref{ex_presheaves-on-sets}(ii).
\end{proof}

\begin{proposition}\label{prop_was-also-cofinal}
Let $\cC,\cD$ be two categories, $F:\cC\to\cD$ a functor, and
$\phi:X\to Y$ a morphism in $\cC$. The following holds :
\begin{enumerate}
\item
If $F$ commutes with pull backs and $\phi$ is a monomorphism,
then $F\phi$ is a monomorphism.
\item
If $F$ commutes with push-outs and $\phi$ is an
epimorphism, then $F\phi$ is an epimorphism.
\end{enumerate}
\end{proposition}
\begin{proof}(i): After replacing $\sU$ by a larger universe,
we may assume that $\cC$ and $\cD$ have small $\Hom$-sets.
Notice then that the tautological cone of $X\times_YX$
amounts to a pair of morphisms
$$
p_1,p_2:X\times_YX\to h_X
\qquad
\text{in $\cC^\wedge$}
$$
such that $h_\phi\circ p_1=h_\phi\circ p_2$ (see remark
\ref{rem_wishful}(iv)), and the diagonal morphism
$\iota_\phi:h_X\to X\times_YX$ (see \eqref{subsec_codiagonal})
is characterized as the unique morphism in $\cC^\wedge$ such
that $p_1\circ\iota_\phi=\one_{h_X}=p_2\circ\iota_\phi$.
Likewise we can characterize $\iota_{F\phi}$. Now,
by assumption there exist unique morphisms $j$ and
$q_i$ for $i=1,2$, fitting into a commutative diagram
\set\begin{equation}\label{eq_middle-arrow}
{\diagram F_!h_X \ar[r]^-{F_!\iota_\phi} \ar[d] &
F_!(X\times_YX) \ar[r]^-{F_!p_i} \ar[d] & F_!h_X \ar[d] \\
h_{FX} \ar[r]^-{j} & FX\times_{FY}FX \ar[r]^-{q_i} & h_{FX}
\enddiagram}
\end{equation}
whose leftmost and rightmost vertical arrows are the
natural identifications of remark \ref{rem_was-cofinal}(v),
and whose central vertical arrow is the isomorphism
given by definition \ref{def_commute-with-lim}(i).
Taking into account proposition \ref{prop_diagonal-monos}(i),
we see that the assertion will follow, once we know
that $j=\iota_{F\phi}$. However, the composition of
the two top horizontal arrows equals $\one_{F_!h_X}$,
so the composition of the bottom horizontal arrows
is $\one_{h_{FX}}$, and thus we come down to showing :

\begin{claim} $q_1$ and $q_2$ are the two morphisms
defining the tautological cone for $FX\times_{FY}FX$.
\end{claim}
\begin{pfclaim} Notice that, for any $Z\in\Ob(\cC)$,
an element of $X\times_YX(Z)$ is the same as the datum
of a pair of morphisms $\alpha,\beta:Z\to X$ such that
$\phi\circ\alpha=\phi\circ\beta$. The natural morphism
$X\times_YX\to F^\wedge(FX\times_{FY}FX)$ of
\eqref{eq_before-adjoint} sends such a pair $(\alpha,\beta)$
to the pair $(F\alpha,F\beta)$. Therefore, the central
vertical arrow of \eqref{eq_middle-arrow} is the map
that assigns to any $W\in\Ob(\cD)$ and any pair
$(\psi:W\to FZ,(\alpha,\beta:Z\to X))$ the pair
$(F\alpha\circ\psi,F\beta\circ\psi:W\to FX)$. Then,
$q_1$ (resp. $q_2$) sends such a pair to $F\alpha\circ\psi$
(resp. to $F\beta\circ\psi$). On the other hand,
$F_!p_1$ (resp. $F_!p_2$) sends $(\psi,(\alpha,\beta))$
to the pair $(\psi,\alpha)$ (resp. to $(\psi,\beta)$).
Now the claim follows by inspecting the explicit
definition of the rightmost vertical arrow of
\eqref{eq_middle-arrow}.
\end{pfclaim}

(ii) follows by dualizing the foregoing argument
(details left to the reader).
\end{proof}

\begin{proposition}\label{prop_conser-nd-left-exact}
{\em(i)}\ \
Let $\cC,\cD$ be two categories, $F:\cC\to\cD$ a left exact
(resp. right exact) functor; suppose that $\cC$ is finitely
complete (resp. finitely cocomplete), and consider the
following conditions :
\begin{enumerate}
\alphaenu
\item
$F$ is conservative (see definition {\em\ref{def_equivalence}(ii)})
\item
$F$ is faithful
\item
$F$ reflects epimorphisms (resp. monomorphisms)
\item
$F$ reflects monomorphisms (resp. epimorphisms).
\end{enumerate}
Then {\em (a)$\Rightarrow$(b)$\Rightarrow$(c)$\Rightarrow$(d)}.

{\em(ii)}\ \
Suppose moreover that the isomorphisms in $\cC$ are the
morphisms that are both monomorphisms and epimorphisms.
Then {\em (c)$\Rightarrow$(a)}.
\end{proposition}
\begin{proof} We prove the assertions in case $f$ is
left exact and $\cC$ is finitely complete; the assertions
for the case where $F$ is right exact and $\cC$ is finitely
cocomplete will follow by considering $F^o:\cC^o\to\cD^o$.

(i.a)$\Rightarrow$(i.b): Let $f;g:A\to B$ be two morphisms
in $\cC$; by assumption, the equalizer of $f$ and $g$ is
representable in $\cC$ by some object $E$, and the universal
cone is the datum of a morphism $h:E\to A$ such that
$f\circ h=g\circ h$. Moreover, $FE$ represents the equalizer
of $Ff$ and $Fg$, and $Fh:FE\to FA$ is a universal cone for
this equalizer, by remark \ref{rem_representable-case}(i).
Now, suppose that $Ff=Fg$; then $Fh$ is an isomorphism,
hence the same holds for $h$, and therefore $f=g$.

(i.b)$\Rightarrow$(i.c): Let $f:A\to B$ and $g,h:B\to C$
be three morphisms in $\cC$ with $g\circ f=h\circ f$, and
suppose that $Ff$ is an epimorphism; since
$Fg\circ Ff=Fh\circ Ff$, it follows that $Fg=Fh$, and then
$g=h$, since $F$ is faithful. This shows that $f$ is an
epimorphism.

(i.c)$\Rightarrow$(i.d): Let $f:A\to B$ be a morphism in
$\cC$ such that $Ff:FA\to FB$ is a monomorphism; by assumption
the fibre product $A\times_BA$ is representable in $\cC$,
and if $p_1,p_2:A\times_BA\to A$ are the natural projections,
then $F(A\times_BA)$ represents $FA\times_{FB}FA$, and
$Fp_1,Fp_2:F(A\times_BA)\to FA$ provide a universal cone
for this fibre product. It follows easily that if
$\iota_f:A\to A\times_BA$ is the diagonal of $f$, then
$F\iota_f:FA\to F(A\times_BA)$ is the diagonal of $Ff$,
hence $F\iota_f$ is an isomorphism in $\cD$, by proposition
\ref{prop_diagonal-monos}(i). By assumption, $\iota_f$ is
then an epimorphism, so $f$ is a monomorphism, again by
proposition \ref{prop_diagonal-monos}(i).

(ii): Let $f:A\to B$ be a morphism such that $Ff$ is
an isomorphism; in particular, $Ff$ is an epimorphism,
so by assumption $f$ is an epimorphism. But $Ff$ is also
a monomorphism, hence the same holds for $f$, because we
know already that (i.c)$\Rightarrow$(i.d). Then $f$ is
an isomorphism, under our assumption on $\cC$.
\end{proof}

\begin{proposition}\label{prop_commute-criteria}
Let $\cC,\cD$ be any two categories, and $f:\cC\to\cD$ any
functor.
\begin{enumerate}
\item
Suppose that
\begin{enumerate}
\item
$f$ commutes with equalizers and with products (resp. and
with finite products)
\item
all small products (resp. all finite products) are representable
in $\cC$.
\end{enumerate}
Then $f$ commutes with all limits of\/ $\cC$ (resp. $f$ is
left exact).
\item
Dually, suppose that
\begin{enumerate}
\item
$f$ commutes with equalizers and with coproducts (resp. and
with finite coproducts)
\item
all small coproducts (resp. all finite coproducts) are
representable in $\cC$.
\end{enumerate}
Then $f$ commutes with all colimits of\/ $\cC$ (resp. $f$ is
right exact).
\end{enumerate}
\end{proposition}
\begin{proof}(i): With the notation of the proof of
proposition \ref{prop_complete-criteria}(i), we have shown
that the rule
$$
(\mu:c_X\Rightarrow F)\mapsto\omega_P^{-1}(\mu*\mathrm{ob})
\qquad
\text{for every $X\in\Ob(\cC)$ and $\mu\in\lim_IF(X)$}
$$
defines an isomorphism $g:\lim_IF\isom E\subset h_P$. Also,
since by assumption $f$ commutes with products (resp. with
finite products), the isomorphisms $\omega_P$ and $\omega_Q$
induce isomorphisms
$$
\gamma_P:h_{fP}\isom\lim_{\Ob(I)}f\circ F\circ\mathrm{ob}
\qquad
\gamma_Q:h_{fQ}\isom
\lim_{\rMorph(I)}f\circ F\circ\mathrm{ob}\circ\mathrm{t}
$$
as explained in remark \ref{rem_representable-case}(i).
We set $E':=\Equal(f\alpha,f\beta)$, and we identify $E'$
as well with a subobject of $h_{fP}$, via the monomorphism
$E'\to h_{fP}$ provided by the tautological cone. Then we
define likewise a morphism $g':\lim_I(f\circ F)\to E'$ by
the rule :
$$
(\mu:c_Y\Rightarrow f\circ F)\mapsto\gamma_P^{-1}(\mu*\mathrm{ob})
\qquad
\text{for every $Y\in\Ob(\cD)$ and $\mu\in\lim_I(f\circ F)(Y)$}.
$$
Consider now the diagram of presheaves on $\cD$ :
\set\begin{equation}\label{eq_products-equals}
{\diagram
f_!\lim_IF \ar[r]^-{f_!g} \ar[d]_d & f_!E \ar[d]^e \\
\lim_If\circ F \ar[r]^-{g'} & E'
\enddiagram}
\end{equation}
whose vertical arrows are the morphisms \eqref{eq_F-commutes-with-f}.

\begin{claim}\label{cl_products-equals}
Diagram \eqref{eq_products-equals} commutes.
\end{claim}
\begin{pfclaim} Let $\tau:c_P\Rightarrow F\circ\mathrm{ob}$ be
the universal cone arising from the isomorphism $\omega_P$;
according to remark \ref{rem_representable-case}(i), for every
$X\in\Ob(\cC)$ and every cone $\mu:c_X\Rightarrow F$ (resp.
every $Y\in\Ob(\cD)$ and every cone
$\mu':c_Y\Rightarrow f\circ F\circ\mathrm{ob}$), the morphism
$g(\mu)$ (resp. $g'(\mu')$) is characterized as the unique one
such that
$$
\tau\odot c_{g(\mu)}=\mu*\mathrm{ob}
\qquad
\text{(resp.\ $(f*\tau)\odot c_{g'(\mu')}=\mu'*\mathrm{ob}$\ )}.
$$
Now, let $(\psi:Y\to fX,\mu:c_X\Rightarrow F)$ be the
representative of any given element $s\in f_!\lim_IF(Y)$.
Then $f_!g(s)$ is represented by the pair $(\psi,g(\mu))$,
and
$$
d(s)=(f*\mu)\odot c_\psi
\qquad
e(f_!g(s))=f(g(\mu))\circ\psi
$$
(remark \ref{rem_was-cofinal}(iii)), and the foregoing
yields the identities
$$
\begin{aligned}
(f*\tau)\odot c_{g'((f*\mu)\circ c_\psi)}
& =((f*\mu)\odot c_\psi)*\mathrm{ob} \\
& =(f*\mu*\mathrm{ob})\odot c_\psi \\
& =(f*(\tau\odot c_{g(\mu)}))\odot c_\psi \\
& =(f*\tau)\odot(f*c_{g(\mu)})\odot c_\psi \\
& =(f*\tau)\odot c_{f(g(\mu))\circ\psi}
\end{aligned}
$$
whence $e(f_!g(s))=g'(d(s))$, as claimed.
\end{pfclaim}

Now, by assumption both $f_!g$ and $e$ are isomorphisms;
in light of claim \ref{cl_products-equals}, it then suffices
to check that $g'$ is an isomorphism. By unwinding the
definitions, we see that the latter assertion amounts to
the following. For every $Y\in\Ob(\cD)$ he rule
$\mu\mapsto\mu*\mathrm{ob}$ yields a bijection from
the set of all cones $\mu:c_Y\Rightarrow f\circ F$ onto
the set of all cones
$\lambda:c_Y\Rightarrow f\circ F\circ\mathrm{ob}$ such that
$\lambda*\mathrm{t}=
((f\circ F)*\mathrm{m})\odot(\lambda*\mathrm{s})$. In turns,
this follows by a direct inspection.
\end{proof}

\begin{proposition}\label{prop_was-get-maddd}
Let $\cA,\cB$ be two categories, $F:\cA\to\cB$ a functor
that admits a left adjoint $G:\cB\to\cA$. The following holds :
\begin{enumerate}
\item
If $\cA$ and $\cB$ are small, we have natural isomorphisms
of functors
$$
F^\wedge\isom G_*
\qquad
F_!\isom G^\wedge.
$$
\item
If $\cA$ and $\cB$ have small $\Hom$-sets, then for every
small category $I$ and every functor $\phi:I\to\cA$, there
is a natural isomorphism of presheaves
$$
G^\wedge(\lim_I\phi)\isom\lim_I\,(F\circ\phi).
$$
\item
$F$ commutes with all limits of $\cA$.
\item
Dually, $G$ commutes with all colimits of $\cB$.
\item
$F$ transforms monomorphisms into monomorphisms, and
$G$ transforms epimorphisms into epimorphisms.
\end{enumerate}
\end{proposition}
\begin{proof}(i): In view of example
\ref{ex_presheaves-on-sets}(ii), the assertion means
that $G^\wedge$ is left adjoint to $F^\wedge$; however
$F^o$ is left adjoint to $G^o$ (remark
\ref{rem_opposite-Fun}(iv)), so it suffices to apply
remark \ref{rem_opposite-Fun}(iii), which shows more
precisely that $\eps^\wedge$ and $\eta^\wedge$ are
respectively the unit and counit of a unique adjunction
for the pair $(G^\wedge,F^\wedge)$.

(ii): We are easily reduced to the case where $\cA$ and
$\cB$ are small. Set $L:=\lim_IF\circ\phi$.
According to remark \ref{rem_wishful}(i), we have
a natural morphism $\omega:\lim_I\phi\to F^\wedge L$,
to which the adjunction exhibited in the proof of (i)
attaches a unique morphism
$\omega^*:G^\wedge(\lim_I\phi)\to L$. We shall show
more precisely that $\omega^*$ is an isomorphism of
presheaves. Indeed, let $B\in\Ob(\cB)$ be any
object; by unwinding the definitions, we see that
$\omega^*_B$ is given by the rule
$$
\mu\mapsto
\mu^*:=((F*\mu)\odot c_{\eta_B}:c_B\Rightarrow F\circ\phi)
\qquad
\text{for every cone $\mu:c_{GB}\Rightarrow\phi$}
$$
(notation of \eqref{subsec_wishful}). We provide an
explicit inverse for the rule $\mu\mapsto\mu^*$, by
setting
$$
\tau^*:=(\eps*\phi)\odot(G*\tau):c_{GB}\Rightarrow\phi
\qquad
\text{for every cone $\tau:c_B\Rightarrow F\circ\phi$}.
$$
Clearly the rule $\tau\mapsto\tau^*$ defines a morphism
of presheaves $L\to G^\wedge(\lim_I\phi)$. Then, for every
such $\tau$ the triangular identities of \eqref{subsec_adj-pair}
give :
$$
\begin{aligned}
(\tau^*)^*
& =(F*((\eps*\phi)\odot(G*\tau)))\odot c_{\eta_B} \\
& =(F*\eps*\phi)\odot(F*G*\tau)\odot c_{\eta_B} \\
& =(F*\eps*\phi)\odot(\eta*F*\phi)\odot\tau \\
& =\tau.
\end{aligned}
$$
Likewise one sees that $(\mu^*)^*=\mu$ for every cone
$\mu:c_{GB}\Rightarrow\phi$, as needed.

(iii): Again, we can assume that $\cA$ and $\cB$ are small.
We have just seen that the adjunction for the
pair $(G^\wedge,F^\wedge)$ exhibited in (i) transforms
the morphism $\omega$ into an isomorphism $\omega^*$.
By example \ref{ex_presheaves-on-sets}(ii), it follows
that any adjunction for the pair $(F_!,F^\wedge)$ will
also transform $\omega$ into an isomorphism, whence
the contention.
(Notice that $F_!$ is -- {\em a priori} -- well defined
only as a functor on $\sV$-presheaves, for some universe
$\sV$ such that $\cA$ and $\cB$ are small $\sV$-categories,
but (i) implies that actually it is -- up to natural
isomorphism -- already defined on $\sU$-presheaves.)

(iv) follows as usual, by considering the opposite
categories and taking into account remark
\ref{rem_opposite-Fun}(iv). Lastly, (v) follows from
(iii),(iv) and proposition \ref{prop_was-also-cofinal}.
\end{proof}

Let $\cA$ be a subcategory of a category $\cB$, and
$\phi:I\to\cA$ a functor from a small category $I$,
whose limit is representable in $\cA$ by an object $A$.
In this situation -- and even in case $\cA$ is a full
subcategory of $\cB$ -- it is not necessarily true that
$A$ represents the limit of the composition $I\to\cB$ of
$\phi$ with the inclusion functor $\cA\to\cB$. We have
nevertheless the following result :

\begin{corollary}\label{cor_mitchell}
Let $\cA,\cB$ be two categories, $F:\cA\to\cB$ a fully
faithful functor admitting a left adjoint $G:\cB\to\cA$,
and $\phi:I\to\cA$ a functor from a small category $I$.
We have :
\begin{enumerate}
\item
If the limit of\/ $F\circ\phi$ is representable by some
$L\in\Ob(\cB)$, then the limit of\/ $\phi$ is representable
by some object $M\in\Ob(\cA)$.
\item
In the situation of\/ {\em (i)}, there exist isomorphisms
$L\isom FM$ in $\cB$ and $M\isom GL$ in $\cA$.
\item
Especially, if $\cB$ is complete, the same holds for $\cA$.
\end{enumerate}
\end{corollary}
\begin{proof} Let $\eta:\one_\cB\Rightarrow F\circ G$ be a unit
and $\eps:G\circ F\Rightarrow\one_\cA$ a counit for the pair $(G,F)$,
and $\tau:c_L\Rightarrow F\circ\phi$ a universal cone; set
$\tau^*:=(\eps*\phi)\odot(G*\tau):c_{GL}\Rightarrow\phi$.
There exists a unique morphism $f:FGL\to L$ in $\cB$ such that
$\tau\odot c_f=F*\tau^*$. We have :
$$
\begin{aligned}
\tau\odot c_{f\circ\eta_L}=\tau\odot c_f\odot c_{\eta_L}=
(F*\tau^*)\odot(\eta*c_L)&\,=
(F*\eps*\phi)\odot(FG*\tau)\odot(\eta*c_L) \\
&\,=(F*\eps*\phi)\odot(\eta*F\phi)\odot\tau=\tau
\end{aligned}
$$
whence $f\circ\eta_L=\one_L$, by the universality of $\tau$.
On the other hand, since $F$ is fully faithful, there exists
a morphism $h:GL\to GL$ in $\cA$ such that
$Fh=\eta_L\circ f:FGL\to FGL$. Now, the adjunction of the pair
$(G,F)$ assigns to $h$ the morphism
$Fh\circ\eta_L=\eta_L\circ f\circ\eta_L=F(\one_{GL})\circ\eta_L$
(see \eqref{subsec_adj-pair}), so we must have $h=\one_{GL}$,
therefore $f$ is an isomorphism. Lastly, since $F$ is
fully faithful, we deduce easily that $GL$ represents
the limit of $F$, whence the corollary.
\end{proof}

\subsection{Special properties of the categories of presheaves}
In this section we gather some results concerning the
structure of categories of presheaves and the exactness
properties of various natural functors between such categories.
We begin with a corollary of proposition \ref{prop_was-get-maddd} :

\begin{corollary}\label{lem_presheaf-in-a-cat}
Let $\cA$ and $\cB$ be two categories. We have :
\begin{enumerate}
\item
If $\cA$ is small and $\cB$ has small $\Hom$-sets, the category
$\bFun(\cA,\cB)$ has small $\Hom$-sets.
\item
Suppose that $\cB$ is complete or cocomplete. Then, for every
$A\in\Ob(\cA)$ the functor
$$
\eps_A:\bFun(\cA,\cB)\to\cB
\qquad
F\mapsto FA
\qquad
(\eta:F\Rightarrow G)\mapsto(\eta_A:FA\to GA)
$$
commutes with all colimits and with all limits.
\item
If $\cB$ is finitely complete (resp. finitely cocomplete)
then $\eps_A$ commutes with finite limits (resp. with
finite colimits), for every $A\in\Ob(\cA)$.
\item
If $\cB$ is complete (resp. cocomplete, resp. finitely
complete, resp. finitely cocomplete), the same holds for
$\bFun(\cA,\cB)$.
\end{enumerate}
\end{corollary}
\begin{proof}(i): If $F,G:\cA\to\cB$ are any two functors,
then $\Hom_{\bFun(\cA,\cB)}(F,G)$ is a subset of the product
$\prod_{A\in\Ob(\cA)}\Hom_\cB(FA,GA)$, hence it is a small set.

(iv): Let $I$ be any small category, and
$$
F:I\to\bFun(\cA,\cB)
\qquad
i\mapsto(F_i:\cA\to\cB)
\qquad
(\phi:i\to j)\mapsto(F_\phi:F_i\Rightarrow F_j)
$$
any functor. For every $A\in\Ob(\cA)$ we obtain a
functor $F_A:I\to\cB$ by the rule
$$
F_A(i):=F_iA
\qquad
F_A(\phi):=(F_\phi)_A:F_iA\to F_jA
$$
for every $i,j\in\Ob(I)$ and every $\phi\in\Hom_I(i,j)$.
If $\cB$ is complete, for every such $A$ we may find
$LA\in\Ob(\cB)$ and a universal cone
$\gamma_A:c_{LA}\Rightarrow F_A$. Next, if $g:A\to A'$ is
any morphism in $\cA$, we obtain a natural transformation
$F_g:F_A\Rightarrow F_{A'}$ by the rule
$$
(F_g)_i:=F_i(g):F_iA\to F_iA'
\qquad
\text{for every $i\in\Ob(I)$}
$$
and the limit of $F_g$ is represented by a morphism
$Lg:LA\to LA'$ in $\cB$ such that
$$
F_g\circ\gamma_A=\gamma_{A'}\circ c_{Lg}.
$$
It is then easily seen that the rules $A\mapsto LA$ and
$g\mapsto Lg$ amount to a well defined functor $\cA\to\cB$
that represents the limit of $F$ : details left to the
reader. Lastly, if $\cB$ is cocomplete (resp. finitely
cocomplete), $\cB^o$ is complete (resp. finitely complete),
so the same holds for $\bFun(\cA^o,\cB^o)$, by the foregoing;
then $\bFun(\cA,\cB)$ is cocomplete (resp. finitely
complete), by remark \ref{rem_opposite-Fun}(i).

(ii): Suppose that $\cB$ is complete (resp. cocomplete);
in this case, the proof of (iv) already shows that $\eps_A$
commutes with all limits (resp. colimits). Hence, it remains
only to check that $\eps_A$ commutes as well with all
colimits (resp. limits). To this aim, we may assume that
$\cA$ and $\cB$ are both small, and we denote by $\one$
the category with one object and one morphism; we let
$i_A:\one\to\cA$ be the (unique) functor that sends the
object of $\one$ to $A$. We have a commutative diagram
of functors
$$
\xymatrix{
\bFun(\cA,\cB) \ar[rr]^-{\bFun(i_A,\cB)} \ar[rrd]_{\eps_A} & &
\bFun(\one,\cB) \ar[d] \\
& & \cB
}$$
whose vertical arrow is an obvious isomorphism of categories.
It then suffices to check that the functor $\bFun(i_A,\cB)$
commutes with all colimits (resp. limits). However, by theorem
\ref{th_Kan-ext}, the functor $\bFun(i_A,\cB)$ admits a right
(resp. left) adjoint, so it suffices to invoke proposition
\ref{prop_was-get-maddd}(ii,iii).

(iii) follows likewise from the proof of (iv).
\end{proof}

\begin{remark}\label{rem_prolyx}
(i)\ \
By inspecting the proof of theorem \ref{th_Kan-ext}, it is
easily seen that, more generally, the functors $\eps_A$ of
corollary \ref{lem_presheaf-in-a-cat} commute with limits
when $\cB$ is a category with representatlbe coproducts.
Dually, if all products of $\cB$ are representable, $\eps_A$
commutes with colimits.

(ii)\ \
In the same vein, suppose moreover that $\Hom_\cA(A,A')$
is a finite set for every $A,A'\in\Ob(\cA)$. Then, again by
inspection of the proof of theorem \ref{th_Kan-ext} we
see that for the functors $\eps_A$ to commute with limits
(resp. colimits) it suffices that all finite coproducts
(resp. finite products) of $\cB$ are representable.

(iii)\ \
Furthermore, the proof of corollary
\ref{lem_presheaf-in-a-cat}(iv) shows more precisely
the following. If $F:I\to\Fun(\cA,\cB)$ is a functor
from a small category $I$, and the limit (resp. colimit)
of $\eps_A\circ F$ is representable in $\cB$ for every
$A\in\Ob(\cA)$, then the limit (resp. colimit) of $F$
is representable in $\Fun(\cA,\cB)$, and $\eps_A$
commutes with this limit (resp. colimit), for every
$A\in\Ob(\cA)$.
\end{remark}

\begin{corollary}\label{cor_pre-misc}
Let\/ $\cC$ be any small category, and $f:F\to G$ a morphism
in $\cC^\wedge$.
\begin{enumerate}
\item
$\cC^\wedge$ has small $\Hom$-sets and is both well-powered
and co-well-powered (see \eqref{subsec_well-power}).
\item
For every $X\in\Ob(\cC)$, the functor
$$
\cC^\wedge\to\Set \quad : \quad F\mapsto FX
\qquad
(f:F\to G)\mapsto(f_X:FX\to GX)
$$
commutes with all limits and all colimits (in other words, the
limits and colimits in $\cC^\wedge$ are computed argumentwise).
\item
$f$ is a monomorphism (resp. an epimorphism) if and only if
the map $f_X:FX\to GX$ is injective (resp. surjective) for
every $X\in\Ob(\cC)$.
\item
$f$ is an isomorphism if and only if it is both a monomorphism
and an epimorphism.
\item
Let $f(F)\in\Ob(\cC^\wedge)$ be the subobject of $F$ with
$f(F)(X):=\Img(f_X:FX\to GX)$ for every $X\in\Ob(\cC)$.
Then $f(F)=\Img(f)=\mathrm{Coim}(f)$ (see example
{\em\ref{ex_pre-misc}}).
\item
The Yoneda embedding $h_\cC:\cC\to\cC^\wedge$ commutes
with all representable limits of\/ $\cC$.
\item
If $\cD$ is any category with small $\Hom$-sets and
$\phi:\cC\to\cD$ any functor, then $\phi^\wedge$ commutes with
all limits and all colimits of\/ $\cD^\wedge$, and $\phi_!$
(resp. $\phi_*$) commutes with all colimits (resp. with all
limits) of\/ $\cC^\wedge$.
\end{enumerate}
\end{corollary}
\begin{proof}(i): The first assertion is a special case
of corollary \ref{lem_presheaf-in-a-cat}(i). Next, to see
that $\cC^\wedge$ is well-powered notice that for every
presheaf $F$ on $\cC$, and every $X\in\Ob(\cC)$, the set
of subsets of $F(X)$ is small, and a subobject of $F$ is
just a compatible system of subsets $F'(X)\subset F(X)$,
for $X$ ranging over the small set of objects of $\cC$.
Likewise one sees that $\cC^\wedge$ is co-well-powered.

(ii) is a special case of corollary \ref{lem_presheaf-in-a-cat}(ii),
but it is also already clear from example \ref{ex_big-complete}(i).

(iii): It is easily seen that if $f_X$ is injective
(resp. surjective) for every $X\in\Ob(\cC)$, then $f$
is a monomorphism (resp. an epimorphism). The converse
follows from (ii) and proposition \ref{prop_was-also-cofinal}.

(iv) is an immediate consequence of (iii).

(v): It follows from (iii) that the induced morphism
$F\to f(F)$ is a quotient of $F$, and the sequence
$F\to f(F)\to G$ is the unique factorization of $f$
(up to unique isomorphism) as a composition of an
epimorphism followed by a monomorphism, whence the
assertion.

(vi) follows directly from example \ref{ex_big-complete}(ii)
and remark  \ref{rem_representable-case}(i).

(vii) is a special case of proposition
\ref{prop_commute-criteria}(iii,iv).
\end{proof}

\begin{example}\label{ex_rep-funct-commute-with-lims}
Let $\cC$ be a category with small $\Hom$-sets, and
$Y\in\Ob(\cC)$ any object.

(i)\ \
Suppose that all the products of small families of objects
of $\cC$ are representable, and consider the functor
$$
Y^\bullet:\Set^o\to\cC^\wedge
$$
which to any ($\sU$-small) set $S$ assigns the product
in $\cC$ of $S$ copies of $Y$, {\em i.e.} the limit of
the functor $\phi^S_Y:S\to\cC$ such that $\phi_Y(s):=Y$
for every $s\in S$ (where $S$ is regarded as a discrete
category), and to any map of sets $f:T\to S$, attaches
the morphism of presheaves $Y^f:=\lim_f\phi_Y$.
The assumption on $\cC$ implies that $Y^\bullet$ is
representable by a functor which we denote as well
$$
Y^\bullet:\Set^o\to\cC
$$
(see remark \ref{rem_represent-morph}(ii)). On the other
hand, we have the functor $h_Y:\cC^o\to\Set$ such that
$h_Y(X^o):=\Hom_\cC(X,Y)$ for every $X\in\Ob(\cC)$, and
notice the natural identifications :
$$
\Hom_{\cC^o}((Y^S)^o,X^o)=\Hom_\cC(X,Y^S)\isom
\Hom_\Set(S,h_Y(X^o))
$$
which say that $h_Y$ is right adjoint to the functor
$(Y^\bullet)^o$. By proposition \ref{prop_was-get-maddd}(iii)
it follows that $h_Y$ commutes with all small limits of $\cC^o$.

(ii)\ \
Likewise, if all the small coproducts of $\cC$ are representable,
we may apply (i) to the functor $h_{Y^o}:\cC\to\Set$. We find
that $h_{Y^o}$ commutes with all small limits of $\cC$.

(iii)\ \
Especially, let $I$, $\cC$ be two small categories,
$F:I\to\cC$ a functor, and $G$ any presheaf on $\cC$.
From (i),(ii), example \ref{ex_big-complete}(i) and
corollary \ref{cor_pre-misc}(vi) we get natural isomorphisms :
\set\begin{equation}\label{eq_colim-and-Hom}
\lim_{I^o}\Hom_{\cC^\wedge}(h^o_\cC\circ F^o,G)\isom
\Hom_{\cC^\wedge}(\colim_I h_\cC\circ F,G)
\end{equation}
\set\begin{equation}\label{eq_lim-and-Hom}
\lim_I\Hom_{\cC^\wedge}(G,h_\cC\circ F)\isom
\Hom_{\cC^\wedge}(G,\lim_IF)
\end{equation}
(more precisely, the limit on the left is represented by
the set on the right, in both cases). Furthermore,
since the tautological cone $\tau$ (resp. the tautological
cocone $\mu$) for $F$ is universal, we get the universal
cone $\Hom_{\cC^\wedge}(\tau,G)$ (resp. $\Hom_{\cC^\wedge}(G,\mu)$)
for the limit appearing in \eqref{eq_colim-and-Hom} (resp.
in \eqref{eq_lim-and-Hom} : see remarks \ref{rem_wishful}(iv))
and \ref{rem_representable-case}(i)).
\end{example}

\sset\subsubsection{}\label{subsec_category-of-elements}
It is an important fact that every presheaf on a small
category is the colimit of a system of representable
presheaves; more generally, for any category $\cC$ and
every $F\in\Ob(\cC^\wedge)$, let us consider the
{\em category of elements of $F$} :
$$
\cFib(F)
$$
whose objects are the pairs $(X,s)$ where $X\in\Ob(\cC)$
and $s\in FX$ (the notation shall be explained in
\eqref{subsec_fibred-cats-II}). The morphisms
$f:(X,s)\to(Y,t)$ in $\cFib(F)$ are the morphisms
$f:X\to Y$ of $\cC$ such that $Ff(t)=s$. We have an obvious
{\em source functor}
$$
\ss_F:\cFib(F)\to\cC
\qquad
(X,s)\mapsto X
\qquad
((X,s)\xrightarrow{f}(Y,t))\mapsto(X\xrightarrow{f}Y).
$$
Notice that if $\cC$ is small, the same holds for
$\cFib(F)$. Also, if $\cC$ has small $\Hom$-sets, Yoneda's
lemma naturally identifies $\cFib(F)$ with the category
$h_\cC\cC/F$ (notation of \eqref{subsec_fibreovercat} and
\eqref{subsec_yoneda}); namely, a pair $(X,s)$ corresponds
to the unique morphism of presheaves
$$
h_{(X,s)}:h_X\to F
\qquad\text{such that}\qquad
h_{(X,s),X}(\one_X)=s.
$$
Moreover, in this case $\ss_F$ is identified with the
source functor $\ss_F:h_\cC\cC/F\to\cC$ of
\eqref{subsec_fibreovercat}. We have a natural cocone
$$
h^F:h_\cC\circ\ss_F\Rightarrow c_F
\qquad
(X,s)\mapsto h_{(X,s)}.
$$

\begin{lemma}\label{lem_lable}
With the notation of \eqref{subsec_category-of-elements},
for every category $\cC$ and every presheaf $F$ on $\cC$,
the cocone $h^F$ is universal.
\end{lemma}
\begin{proof} Let $G$ be any presheaf on $\cC$. A natural
transformation $\tau:h_\cC\circ\ss_F\Rightarrow c_G$ is a
rule that assigns to every $(X,s)\in\Ob(h_\cC\cC/F)$ a
morphism of presheaves $h_X\to G$; the latter is naturally
identified with an element $\tau_{(X,s)}\in GX$. Thus, $\tau$
induces a map $\tau_{(X,\bullet)}:FX\to GX$ for every
$X\in\Ob(\cC)$, by the rule : $(X,s)\mapsto\tau_{(X,s)}$.
Moreover, the naturality of $\tau$ implies that
$$
Gf(\tau_{(Y,s)})=\tau_{(X,s)}
\qquad
\text{for every morphism $f:(X,s)\to (Y,t)$ in $h_\cC\cC/F$}
$$
(details left to the reader); {\em i.e.} the rule
$X\mapsto\tau_{(X,\bullet)}$ for every $X\in\Ob(\cC)$ yields a
morphism of presheaves $h^\tau:F\to G$, and it is easily seen
that $\tau=c_{h^\tau}\odot h^F$, whence the contention.
\end{proof}

\begin{example}\label{ex_lim-yoneda}
Let $\cC$ be a small category, and $\bone_\cC$ a final
object of $\cC^\wedge$  (see example \ref{ex_equalizers}(v));
since the source functor $\ss_{\bone_\cC}:h_\cC\cC/\bone_\cC\to\cC$
is an isomorphism of categories, we deduce from lemma
\ref{lem_lable} a natural isomorphism
$$
\colim_\cC h_\cC\isom\bone_\cC
\qquad
\text{in $\cC^\wedge$}.
$$
\end{example}

\sset\subsubsection{}\label{subsec_was-ex_localize-topoi}
We conclude this section with a few observations concerning
presheaves on the slice categories of \eqref{subsec_slice-cat}. 
Thus, let $\cC$ be a category, $X$ any object of $\cC$, denote by
$\ss_X:\cC\!/\!X\to\cC$ the source functor of \eqref{subsec_slice-cat},
and suppose that $\cC$ is $\sV$-small, for some universe $\sV$.
To begin with, a direct inspection of \eqref{eq_describe_f_lower}
yields a natural isomorphism :
\set\begin{equation}\label{eq_can-be-used}
(\ss_X)_{\sV!}F(Y)\isom
\{(\phi,a)~|~\phi\in\Hom_\cC(Y,X),\ a\in F(\phi)\}
\end{equation}
for every $\sV$-presheaf $F$ on $\cC\!/\!X$ and every $Y\in\Ob(\cC)$.
Indeed, from {\em loc.cit.} we see that every element of
$(\ss_X)_{\sV!}F(Y)$ is the equivalence class of a datum
$$
(\psi:Y\to Z,\alpha:Z\to X,s\in F(\alpha))
$$
and such a datum is equivalent to
$(\one_Y,\alpha\circ\psi:Y\to X,F(\psi)(s)\in F(\alpha\circ\psi))$.
If $\psi:Z\to Y$ is any morphism of $\cC$, the corresponding map
$(\ss_X)_{\sV!}F(Y)\to(\ss_X)_{\sV!}F(Z)$ is given by the rule :
\set\begin{equation}\label{eq_specify-map}
(\phi,a)\mapsto(\phi\circ\psi,F(\psi)(a))
\qquad
\text{for every $\phi:Y\to X$ and $a\in F(\phi)$}.
\end{equation}
Moreover, under this natural identification the adjoint
pair $((\ss_X)_{\sV§},\ss^\wedge_X)$ assigns to any
morphism $t:F\to\ss^\wedge_XG$ of presheaves on $\cC\!/\!X$
the morphism $t^*:(\ss_X)_{\sV!}F\to G$ in $\cC^\wedge$
given by the rule :
$$
(\phi,a)\mapsto t(\phi)(a)
\qquad
\text{for every $Y\in\Ob(\cC)$, $\phi:Y\to X$ and $a\in F(\phi)$}.
$$

\begin{proposition}\label{prop_in-the-same-vein}
Let $\cC$ be any category and $X$ any object of\/ $\cC$.
\begin{enumerate}
\item
The source functor\/ $\ss_X:\cC\!/\!X\to\cC$ commutes with
all connected limits and all representable colimits of\/
$\cC\!/\!X$.
\item
Dually, the target functor $\st_X:X/\cC\to\cC$ commutes
with all connected colimits and all representable limits
of $X\!/\cC$.
\item
If all the binary products of\/ $\cC$ are representable,
$\ss_X$ admits a right adjoint, and therefore it commutes
with all colimits of\/ $\cC\!/\!X$.
\item
Dually, if all the binary coproducts of\/ $\cC$ are
representable, $\st_X$ admits a left adjoint, and therefore
it commutes with all limits of\/ $X/\cC$.
\item
If\/ $\cC$ is complete (resp. cocomplete, resp. finitely
complete, resp. finitely cocomplete), then the same holds
for both\/ $\cC\!/\!X$ and\/ $X/\cC$.
\item
Suppose that\/ $\cC$ has small $\Hom$-sets. Then the
following holds :
\begin{enumerate}
\item
The functor $\ss_X^\wedge$ admits a left adjoint
$\ss_{X!}:(\cC\!/\!X)^\wedge\to\cC^\wedge$.
\item
Moreover $\ss_{X!}$ factors as the composition
of an equivalence of categories
$$
e_X:(\cC\!/\!X)^\wedge\to\cC^\wedge\!/h_X
$$
and the source functor $\ss_{h_X}:\cC^\wedge\!/h_X\to\cC^\wedge$
(notation of \eqref{subsec_yoneda}).
\item
$\ss_{X!}$ commutes with fibre products and preserves
monomorphisms.
\end{enumerate}
\end{enumerate}
\end{proposition}
\begin{proof}(v): Suppose that $\cC$ is complete (resp.
finitely complete); to show that the same holds for
$\cC\!/\!X$, we may apply the criterion of proposition
\ref{prop_complete-criteria}(i). Indeed, suppose that
$p_\bullet:=(p_i:Y_i\to X~|~i\in I)$ is any family
(resp. any finite family) of objects of $\cC\!/\!X$,
indexed by a small set $I$, which we regard as a discrete
category. We consider the category $I_X$ whose set of
objects is the disjoint union of $I$ and $\{X\}$, and with
$$
\rMorph(I_X)=
\rMorph(I)\cup\{(X,X,\one_X)\}\cup\{(i,X,p_i)~|~i\in I\}
$$
so that $X$ is the unique final object of $I_X$. Next, we
define a functor $F:I_X\to\cC$, by letting
$$
FX:=X
\qquad\text{and}\qquad
Fi:=Y_i
\quad
Fp_i:=p_i
\qquad\text{for every $i\in I$}.
$$
Any universal cone $\tau:c_L\Rightarrow F$ yields a morphism
$\tau_X:L\to X$, and it is easily seen that the object $\tau_X$
of $\cC\!/\!X$ represents the product of the family $p_\bullet$.
Lastly, let $p_i:Y_i\to X$ ($i=1,2$) be two objects of $\cC\!/\!X$,
and $f,g:p_1\to p_2$ two morphisms in $\cC\!/\!X$. Denote by
$E$ any object of $\cC$ representing the equalizer of
$\ss_X(f),\ss_X(g):Y_1\to Y_2$ in $\cC$, where
$\ss_X:\cC\!/\!X\to\cC$ is the source functor. The
corresponding universal cone yields a morphism $h:E\to Y_1$
such that $\ss_X(f)\circ h=\ss_X(g)\circ h$, and it is
easily seen that the object $p_1\circ h:E\to X$ of $\cC\!/\!X$
represents the equalizer of $f$ and $g$ (details left to the
reader).

Next, suppose that $\cC$ is cocomplete (resp. finitely
cocomplete), and let $F:I\to\cC\!/\!X$ be a functor from
a small (resp. finite) category $I$. Let also $C$ be any
object of $\cC$ representing the colimit of $\ss_X\circ F$;
the functor $F$ can be regarded as a cocone
$\ss_X\circ F\Rightarrow c_X$, which corresponds to a
morphism $\tau:C\to X$, and it is easily seen that
$\tau\in\Ob(\cC\!/\!X)$ represents the colimit of $F$.

The assertions for $X/\cC$ follows as usual, by considering
the opposite categories.

(vi.a): The discussion of \eqref{subsec_was-ex_localize-topoi}
yields a left adjoint $(\ss_X)_{\sV!}$ to $(\ss_X)^\wedge_\sV$, for
any universe $\sV$ such that $\cC$ has $\sV$-small $\Hom$-sets.
However, by inspecting \eqref{eq_can-be-used}, we see that
if $\cC$ has small $\Hom$-sets and $F$ is a $\sU$-presheaf
on $\cC\!/\!X$, then $(\ss_X)_{\sV!}F$ is a $\sU$-presheaf.
Since the inclusion $(\cC\!/\!X)^\wedge_\sU\to(\cC\!/\!X)^\wedge_\sV$
is fully faithful, we deduce that the target of \eqref{eq_can-be-used}
can be used to define the sought left adjoint
$(\ss_X)_!:(\cC\!/\!X)^\wedge_\sU\to\cC^\wedge_\sU$ to
$(\ss_X)^\wedge_\sU$.

(vi.b): Notice that the presheaf $\bone_{\cC/X}:=h_{\one_X}$ is a final
object of $(\cC/X)^\wedge$ : indeed, for every object $f:Y\to X$
of $\cC/X$ we have $\bone_{\cC/X}(f)=\{f/X\}$ (notation of
\eqref{subsec_slice-cat}). Then remark \ref{rem_was-cofinal}(v)
yields a natural isomorphism :
$$
\ss_{X!}(\bone_{\cC/X})\isom h_X
$$
that identifies the unit of adjunction with the morphism of
presheaves :
\set\begin{equation}\label{eq_bone-unit}
\bone_{\cC/X}\to\ss^\wedge_X(h_X)
\qquad\text{such that}\qquad
f/X\mapsto f
\qquad\text{for every $f\in\Ob(\cC/X)$}.
\end{equation}
Now, every presheaf $F$ on $\cC/X$ admits a unique morphism
$F\to\bone_{\cC/X}$, whence a natural morphism $\ss_{X!}F\to h_X$,
so we see already that $\ss_{X!}$ factors through $\ss_{h_X}$
and a functor $e_X$ as required; more precisely,
\eqref{eq_specify-map}, says that the morphism
$e_X(F):\ss_{X!}F\to h_X$ is given by the rule :
$(f,a)\mapsto f$ for every $Y\in\Ob(\cC)$ and every
$(f,a)\in\ss_{X!}F(Y)$. To check that $e_X$ is an equivalence,
we exhibit an explicit quasi-inverse : namely, to every
object $\phi:F\to h_X$ of $\cC^\wedge/h_X$ we assign the
presheaf
$$
e'_X(\phi):=\ss^\wedge_X(F)\times_{\ss^\wedge_X(h_X)}\bone_{\cC/X}
$$
defined via \eqref{eq_bone-unit} and the morphism
$\ss^\wedge_X(\phi):\ss^\wedge_X(F)\to\ss^\wedge_X(h_X)$. Now,
let $F$ be any presheaf on $\cC/X$, and $f:Y\to X$ any
object of $\cC/X$; by unwinding the definitions, we see
that
$$
(e'_X\circ e_X(F))(f)=\ss_{X!}F(Y)\times_{h_X(Y)}\{f/X\}
$$
which is the set of all pairs $(f,a)$ with $a\in Ff$.
Thus we get an obvious isomorphism of presheaves
$e'_X\circ e_X(F)\isom F$, natural with respect to $F$.
Likewise, for every object $\phi:G\to h_X$ of
$\cC^\wedge/h_X$ and every $Y\in\Ob(\cC)$ we get the map
$$
(e_X\circ e'_X(\phi))_Y:
\ss_{X!}(\ss^\wedge_X(F)\times_{\ss^\wedge_X(h_X)}\bone_{\cC/X})(Y)\to h_X(Y)
$$
which is described as follows. First,
$\ss_{X!}(\ss^\wedge_X(F)\times_{\ss^\wedge_X(h_X)}\bone_{\cC/X})(Y)$
is the set of all pairs $(f,a)$ with $f:Y\to X$ a morphism
of $\cC$, and $a\in GY\times_{h_X(Y)}\{f/X\}$; {\em i.e.}
the pairs $(f,a)$ with $\phi_Y(a)=f$. Then, the map
$(e_X\circ e'_X(\phi))_Y$ is given by the rule : $(f,a)\mapsto f$
for every such pair $(f,a)$. Again, we deduce an obvious
isomorphism $e_X\circ e'_X(\phi)\isom\phi$ of $h_X$-objects
of $\cC^\wedge$, natural with respect to $\phi$, as required.

(vi.c) follows from the explicit description of $\ss_{X!}$
given in \eqref{subsec_was-ex_localize-topoi}, together with
corollary \ref{cor_pre-misc}(ii,iii).

(i): The assertion for representable colimits follows from
the proof of (v) (together with remark
\ref{rem_representable-case}(i)). Next, let $F:I\to\cC\!/\!X$
be any functor from a connected small category $I$; taking
into account examples \ref{ex_presheaves-on-sets}(ii) and
\eqref{subsec_was-ex_localize-topoi} it suffices to show
that the natural morphism
$$
((\ss_X)_{\sV!}\lim_IF)(Z)\to(\lim_I\ss_X\circ F)(Z)
\qquad
(\phi:Z\to X,\tau:c_\phi\Rightarrow F)\mapsto
(\ss_X*\tau:c_Z\Rightarrow\ss_X\circ F)
$$
is an isomorphism for every $Z\in\Ob(\cC)$ and any sufficiently
large universe $\sV$. Now, let $\tau:c_Z\Rightarrow\ss_X\circ F$
be any cone (with vertex an arbitrary $Z\in\Ob(\cC)$); we remark :

\begin{claim}\label{cl_lift-cones}
For every $i,j\in\Ob(I)$ we have
$Fi\circ\tau_i=Fj\circ\tau_j:Z\to X$.
\end{claim}
\begin{pfclaim} Since $I$ is connected, a simple induction
argument reduces to the case where there exists $k\in\Ob(I)$
with morphisms $a:k\to i$ and $b:k\to j$. In this case, we
may compute directly :
$Fi\circ\tau_i=Fi\circ\ss_X(Fa)\circ\tau_k=F_k\circ\tau_k=
Fj\circ\ss_X(Fb)\circ\tau_k=Fj\circ\tau_j$.
\end{pfclaim}

Hence, pick any $i\in\Ob(I)$ (recall that this set is non-empty),
and set $\phi:=Fi\circ\tau_i$; claim \ref{cl_lift-cones} says
that $\tau$ lifts uniquely to a cone $\tau_{/X}:c_\phi\Rightarrow F$,
such that $\ss_X*\tau_{/X}=\tau$, whence the contention.

(iv): Pick a universe $\sV$ such that $\cC$ has $\sV$-small
$\Hom$-sets. For every $Y\in\Ob(\cC)$, the functor
$$
\cC\to\cC^\wedge_\sV
\qquad
Y\mapsto X\amalg Y
$$
is then representable by a functor $\cC\to\cC$, and
we abuse notation, by denoting both functors in the
same way. Moreover, the universal cocone for the
coproduct consists of a pair of morphisms
$X\xrightarrow{\ i_{X,Y}\ }X\amalg Y\xleftarrow{\ i'_{X,Y}\ }Y$
whence a functor
$$
i_X:\cC\to X/\cC
\quad : \quad
Y\mapsto i_{X,Y}
\qquad
\text{for every $Y\in\Ob(\cC)$}.
$$
It is easily seen that $i_X$ is left adjoint to $\st_X$
(details left to the reader). Thus, the assertion follows
from proposition \ref{prop_was-get-maddd}(i,iii).

(iii): Notice that the source functor
$\ss_X:\cC\!/\!X\to\cC$ equals $\st^o_{X^o}$ and the
coproducts of $\cC^o$ are representable; then the
assertion follows from (iv) and remark
\ref{rem_opposite-Fun}(iv). The same argument shows
that (i) implies (ii).
\end{proof}

\begin{remark}\label{rem_in-the-same-vein}
Let $\cC$ and $X$ be as in proposition
\ref{prop_in-the-same-vein}(vi).

(i)\ \
The functor $\ss_{X!}$ does not generally preserve final
objects, hence it is not generally exact.

(ii)\ \
The quasi-inverse to $e_X$ constructed in the proof
of proposition \ref{prop_in-the-same-vein}(vi.b), can
be described more compactly as the functor that assigns
to any object $g:F\to h_X$ of $\cC^\wedge\!/h_X$ the
presheaf given by the rule :
$$
(Y\xrightarrow{\ \phi\ } X)\mapsto
\Hom_{\cC^\wedge\!/h_X}((h_Y\xrightarrow{\ h_\phi\ } h_X),g).
$$
\end{remark}

\begin{example}\label{ex_universal-col}
(i)\ \
Let $\cC$ be a category whose binary products are
representable, and pick a universe $\sV$ such that
$\cC$ has $\sV$-small $\Hom$-sets. For every
$X\in\Ob(\cC)$, we have a functor
$$
p_X:\cC\to\cC^\wedge_\sV
\qquad
Y\mapsto X\times Y
$$
and the assumption on $\cC$ implies that $p_X$ is
representable by a functor $\pi_X:\cC\to\cC$ (remark
\ref{rem_represent-morph}(ii)). Let now $F:I\to\cC$
be a functor from a small category $I$; we say that
{\em the colimit of $F$ commutes with products}, if
$\pi_X$ commutes with the colimit of $F$, for every
$X\in\Ob(\cC)$.

(ii)\ \
A more restricted class of colimits plays a role
in applications. Namely, suppose now that all fibre
products of $\cC$ are representable. Let $F$ be
as in (i), and consider any cone $\tau:F\Rightarrow c_X$,
for any $X\in\Ob(\cC)$. Notice that $\tau$ is the
same as the datum of a functor $F_\tau:I\to\cC\!/\!X$
which {\em lifts} $F$, {\em i.e.} such that
$\ss_X\circ F_\tau=F$, where $\ss_X:\cC\!/\!X\to\cC$
is the source functor. Moreover, the assumption on
$\cC$ implies that the products are representable
in $\cC\!/\!X$, for every such $X$. In this situation,
we shall then say that the colimit of $F$ is
{\em universal}, if the colimit of $F_\tau$ commutes
with products, for every $X$ and $\tau$ as above.

(iii)\ \
Suppose moreover that the colimit of $F$ is representable
by an object $C^o$ of $\cC^o$, and fix a universal cocone
$\mu:F\Rightarrow c_C$. Then any $\tau$ as in (ii) corresponds
to a unique morphism $f:C\to X$, and $\mu$ lifts to a
universal cocone $\mu_\tau:F_\tau\Rightarrow c_f$. For
every morphism $g:Y\to X$ in $\cC$ we may consider the
functor
$$
F\times_{(\tau,g)}Y:=\ss_X\circ\pi_g\circ F_\tau:I\to\cC
$$
(that is, $F\times_{(\tau,g)}Y(i)$ represents $Fi\times_{(\tau_i,g)}Y$
for every $i\in\Ob(I)$). Likewise, we abuse notation and
write $C\times_{(f,g)}Y$ to denote the representative
$\ss_X\circ\pi_g(f)$ for this fibre product; taking into
account proposition \ref{prop_in-the-same-vein}(v), we see
that the colimit of $F$ is universal if and only if the
cocone
$$
(\ss_X\circ\pi_g)*\mu_\tau:F\times_{(\tau,g)}Y\Rightarrow
c_{C\times_{(f,g)}Y}
$$
is universal for every such $X$, $\tau$ and $g$. 

(iv)\ \
For instance, from the explicit descriptions of example
\ref{ex_complete-cats}(i) one may deduce that all colimits
in $\Set$ are universal. In view of example \ref{ex_big-complete}(i)
and corollary \ref{cor_pre-misc}(ii), it follows easily that all
colimits in $\cB^\wedge$ are universal, for any small category $\cB$.
\end{example}

\sset\subsubsection{}\label{subsec_bernard}
Let $F:\cA\to\cB$ be a functor from a small category $\cA$
to a category $\cB$ with small $\Hom$-sets. By remark
\ref{rem_was-cofinal}(ii), for every object $A\in\Ob(\cA)$
we have a natural isomorphism $h_{FA}\isom F_!A$; combining
with example \ref{ex_comma-adjunction}(i), we deduce that
the adjoint pair of functors $F^\wedge:\cB^\wedge\to\cA^\wedge$
and $F_!:\cA^\wedge\to\cB^\wedge$ induces an adjoint pair :
$$
(F^\wedge)_{|h_A}:\cB^\wedge/h_{FA}\isom\cB^\wedge/F_!h_A\to\cA^\wedge/h_A
\qquad
(F_!)_{|h_A}:\cA^\wedge/h_A\to\cB^\wedge/F_!h_A\isom\cB^\wedge/h_{FA}.
$$
For every category $I$, let $I_\circ$ be the category such that
\begin{itemize}
\item
$\Ob(I_\circ)$ is the disjoint union of $\Ob(I)$ with a set
$\{i_\circ\}$
\item
$i_\circ$ is the unique final object of $I_\circ$
\item
$I$ is a full subcategory of $I_\circ$ (more precisely, the
full subcategory of $I_\circ$ whose set of objects is $\Ob(I)$
coincides with $I$).
\end{itemize}
It is clear that these conditions determine uniquely $I_\circ$.
Moreover, for every category $\cC$ and every $C\in\Ob(\cC)$,
every functor $\phi:I\to\cC/C$ extends uniquely to a well
defined functor $\phi_\circ:I_\circ\to\cC/C$ such that
$\phi_\circ(i_\circ)=\one_C$.

\begin{proposition}\label{prop_done-right}
{\em (i)}\ \
With the notation of \eqref{subsec_bernard}, for every
$A\in\Ob(\cA)$ we have {\em essentially commutative} diagrams :
$$
\xymatrix@C+20pt{
(\cB/FA)^\wedge \ar[d]_{e_{FA}} \ar[r]^-{(F_{|A})^\wedge} &
(\cA/A)^\wedge \ar[d]^{e_A} &
(\cA/A)^\wedge \ar[r]^-{(F_{|A})_!} \ar[d]_{e_A} &
(\cB/FA)^\wedge \ar[d]^{e_{FA}} \\
\cB^\wedge/h_{FA} \ar[r]^-{(F^\wedge)_{|h_A}} & \cA^\wedge/h_A &
\cA^\wedge/h_A \ar[r]^-{(F_!)_{|h_A}} & \cB^\wedge/h_{FA}
}$$
where $e_A$ and $e_{FA}$ are the equivalences of proposition
{\em\ref{prop_in-the-same-vein}(vi.b)}.

{\em (ii)}\ \
Let moreover $\phi:I\to\cA/A$ be any functor. Then $F_{|A}$
commutes with the limit of $\phi$ if and only if $F$ commutes
with the limit of $\ss_A\circ\phi_\circ:I_\circ\to\cA$.
\end{proposition}
\begin{proof}(i): The assertion means that there exists an
isomorphism of functors
$e_A\circ(F_{|A})^\wedge\isom(F^\wedge)_{|h_A}\circ e_{FA}$, and
likewise for the second diagram. Let $\phi:K\to h_{FA}$ be
any object of $\cB^\wedge/h_{FA}$; with $e'_A$ and $e'_{FA}$ as
in the proof of proposition \ref{prop_in-the-same-vein}(vi.b),
we get :
$$
\begin{aligned}
(F_{|A})^\wedge\circ e'_{FA}(\phi)
&\,=(F_{|A})^\wedge(\ss_{FA}^\wedge K\times_{\ss^\wedge_{FA}h_{FA}}\bone_{\cB/FA}) \\
&\,\isom(\ss_{FA}\circ F_{|A})^\wedge K
\times_{(\ss_{FA}\circ F_{|A})^\wedge h_{FA}}\bone_{\cA/A} \\
&\,\isom\ss^\wedge_A\circ F^\wedge K
\times_{\ss^\wedge_A\circ F^\wedge h_{FA}}\bone_{\cA/A} \\
&\,\isom e'_A(u^\wedge K\times_{F^\wedge h_{FA}}h_A\to h_A) \\
&\,\isom e'_A\circ(F^\wedge)_{|h_A}(\phi)
\end{aligned}
$$
whence the essential commutativity of the left square diagram.
Since $e_A$ and $e_{FA}$ are equivalences, and since $(F_!)_{|h_A}$
is left adjoint to $(F^\wedge)_{|h_A}$, we deduce the essential
commutativity for the right square diagram as well.

(ii): Let $j:I\to I_\circ$ be the inclusion functor; notice
first that for every category $\cC$, every $C\in\Ob(\cC)$
and every functor $\psi:I\to\cC$ we have $\psi=\psi_\circ j$,
and the morphism of presheaves
$$
\lim_j\one_\cC:\lim_{I_\circ}\psi_\circ\to\lim_I\psi
$$
is an isomorphism (notation of remark \ref{rem_wishful}(i)).
Moreover, it is clear that
$F_{|A}\circ\phi_\circ=(F_{|A}\circ\phi)_\circ$. There follows
a commutative diagram of presheaves :
$$
\xymatrix{ (F_{|A})_!(\lim_{I_\circ}\phi_\circ) \ar[r] \ar[d] &
\lim_{I_\circ}F_{|A}\circ\phi_\circ \ar[d] \\
(F_{|A})_!(\lim_I\phi) \ar[r] & \lim_IF_{|A}\circ\phi.
}$$
whose vertical arrows are isomorphisms, and whose horizontal
arrows are the morphisms of definition
\ref{def_commute-with-lim}(i). Especially, $F_{|A}$ commutes
with the limit of $\phi$ if and only if it commutes with
the limit of $\phi_\circ$. We are thus reduced to checking
the following more general :

\begin{claim} Let $J$ be a connected category, and
$\psi:J\to\cA/A$ a functor. Then $F_{|A}$ commutes with
the limit of $\psi$ if and only if $F$ commutes with
the limit of $\ss_A\circ\psi$. 
\end{claim}
\begin{pfclaim}[] Recall first that
$\ss_{FA!}=\ss_{h_{FA}}\circ e_{FA}$. Since $\ss_{h_{FA}}$ is
obviously conservative, and $e_{FA}$ is an equivalence,
it follows that $\ss_{FA!}$ is conservative. Recall moreover
that $\ss_A$ and $\ss_{FA}$ commute with every connected limit
(proposition \ref{prop_in-the-same-vein}(i)); summing up, we
see that $F_{|A}$ commutes with the limit of $\psi$ if and only
if the morphism
$$
\ss_{FA!}\circ(F_{|A})_!(\lim_J\psi)\xrightarrow{\ss_{FA!}(\omega)}
\ss_{FA!}(\lim_JF_{|A}\circ\psi)\xrightarrow{\omega'}
\lim_J\ss_{FA}\circ F_{|A}\circ\psi
$$
is an isomorphism, where
$\omega:(F_{|A})_!(\lim_J\psi)\to\lim_JF_{|A}\circ\psi$
and $\omega'$ are the morphisms of definition
\ref{def_commute-with-lim}(i). But this composition
is, up to isomorphism, the same as the morphism
\set\begin{equation}\label{eq_toil-and-trouble}
(F\circ\ss_A)_!(\lim_J\psi)\to\lim_JF\circ\ss_A\circ\psi
\end{equation}
of definition \ref{def_commute-with-lim}(i), relative to
the functor $\ss_{FA}\circ F_{|A}=F\circ\ss_A$ and the
limit of $\psi$ (see the proof of proposition
\ref{prop_double-double}(i)). By the same token,
\eqref{eq_toil-and-trouble} is also the composition
$\omega'''\circ F_!(\omega'')$ of the same type, where
$\omega'':\ss_{A!}(\lim_J\psi)\to\lim_J\ss_A\circ\psi$
is an isomorphism (since $J$ is connected), and
$\omega''':F_!(\lim_J\psi)\to\lim_JF\circ\psi$ is an
isomorphism if and only if $F$ commutes with $J$.
\end{pfclaim}
\end{proof}

\subsection{Final and cofinal functors}
Let $I$ be a small category, $\cC$ a category, and $F:I\to\cC$
any functor. For the computation of the limit or colimit of $F$,
it may sometimes be desirable to replace the indexing category
$I$ by simpler ones. That is, we would like to be able to detect
whether a given functor $\phi:J\to I$ induces an isomorphism
from the colimit of $F$ to that of $F\circ\phi$, and if possible,
to construct useful functors of this type, to aid with the
calculation of limits or colimits. Concerning the first aim, one
has a general criterion, for which we shall need the following :

\begin{definition}\label{def_MacL-cofinal}
Let $I$, $J$ be any two categories, and $\phi:J\to I$
any functor. 
\begin{enumerate}
\item
We say that $\phi$ is {\em cofinal} if the category
$i/\phi J$ is connected, for every $i\in\Ob(I)$
(see definition \ref{def_filtered-cols}(ii)).
\item
We say that $\phi$ is {\em final} if $\phi^o:J^o\to I^o$
is cofinal.
\item
If $\phi$ is the inclusion functor of a subcategory
$J$ of $I$, and $\phi$ is cofinal (resp. final) we also
say that $J$ is {\em cofinal in $I$} (resp. {\em final in $I$}).
\item
We say that the category $I$ is {\em finally small} (resp.
{\em cofinally small}), if there exists a final (resp.
cofinal) functor $\phi:J\to I$ with $J$ a small category.
(See corollary \ref{cor_char-finally-small}.)
\end{enumerate}
\end{definition}

\begin{proposition}\label{prop_MacL-cofinal}
Let $\phi:J\to I$ be any functor between small
categories. Then the following conditions are
equivalent :
\begin{enumerate}
\alphaenu
\item
For every category $\cC$ with small $\Hom$-sets, and
every functor $F:I\to\cC$, the induced morphism of
presheaves on $I^o$
\set\begin{equation}\label{eq_cofinal}
\colim_\phi\one_\cC:\colim_IF\to\colim_JF\circ\phi
\end{equation}
is an isomorphism (see remark {\em\ref{rem_wishful}(i)}).
\item
The functor $\phi$ is cofinal.
\end{enumerate}
\end{proposition}
\begin{proof}(a)$\Rightarrow$(b): We fix $i\in\Ob(I)$
and we apply the definition to the functor
$$
h_{i^o}:I\to\Set
\qquad
i'\mapsto\Hom_I(i,i')
$$
whose colimit is a presheaf on the category $\Set^o$.
Then the assertion is an immediate consequence of the
following :

\begin{claim} For every $i\in\Ob(I)$ we have :
\begin{enumerate}
\item
The colimit of $h_{i^o}$ is representable by a set
with one element.
\item
The colimit of $h_{i^o}\circ\phi$ is representable by
$\pi_0(i/\phi J)$ (notation of remark \ref{rem_cofinal}(ii)).
\end{enumerate}
\end{claim}
\begin{pfclaim} Notice that (i) is the special case of
(ii) for $\phi=\one_I$, since the category $i/I$ admits
an initial object, and hence is connected.

(ii): Let $h_I:I\to I^\wedge$ be the Yoneda embedding and
choose a representative $L\in\Ob(I^\wedge)$ for the colimit
of the functor $h_I\circ\phi$. By corollary
\ref{cor_pre-misc}(ii) and remark
\ref{rem_representable-case}(i), the set $L(i)$ represents
the colimit of $h_{i^o}\circ\phi$. On the other hand, from
remark \ref{rem_was-cofinal}(ii) we see that $L$ represents
as well the colimit of $\phi_!\circ h_J$. Moreover, recall
that the colimit of $h_J$ is represented by the final object
$\bone_J$ of $J^\wedge$ (example \ref{ex_lim-yoneda}); by
corollary \ref{cor_pre-misc}(vii) and remark
\ref{rem_representable-case}(i) we deduce that $L$ is
isomorphic to $\phi_!(\bone_J)$. Then the assertion follows
from example \ref{ex_was-lim-yoneda}.
\end{pfclaim}

(b)$\Rightarrow$(a): Let $F:I\to\cC$ be any functor;
recall that the colimit of $F$ is the functor
$L:\cC\to\Set$ that assigns to any $X\in\Ob(\cC)$ the
set $L(X)$ of all cocones $F\Rightarrow c_X$, {\em i.e.}
all the compatible systems of morphisms
$(\tau_i:Fi\to X~|~i\in\Ob(I))$ in $\cC$, such that
$$
\tau_{i'}\circ Fg=\tau_i
\qquad
\text{for every morphism $g:i\to i'$ in $I$}.
$$
Likewise, the colimit of $F\circ\phi$ is the functor
$L':\cC\to\Set$ that assigns to any $X\in\Ob(\cC)$ the
set $L'(X)$ of all compatible systems
$(\mu_j:F\phi(j)\to X~|~j\in\Ob(J))$ in $\cC$, such that
$$
\mu_{j'}\circ F\phi(h)=\mu_j
\qquad
\text{for every morphism $h:j\to j'$ in $J$}.
$$
Under these identifications, \eqref{eq_cofinal} is
the map that assigns to a given $\tau$ as
above, the compatible system $\tau*\phi$,
such that $(\tau*\phi)_j:=\tau_{\phi(j)}$ for every
$j\in\Ob(J)$. We have to check that the rule
$\tau\mapsto\tau*\phi$ is a bijection
$\beta_X:L(X)\isom L'(X)$. However, say that
$\tau*\phi=\tau'*\phi$ for given $\tau,\tau'\in L(X)$,
and consider any $i\in\Ob(I)$; since $i/\phi J$ is
non-empty, we may find $j\in\Ob(J)$ and a morphism
$g:i\to\phi(j)$ in $I$, whence
$$
\tau_i=\tau_{\phi(j)}\circ Fg=\tau'_{\phi(j)}\circ Fg=\tau'_i
$$
hence $\beta_X$ is injective. Next, let $\mu\in L'(X)$
be any compatible system; we consider the map
$$
\tau^*:\Ob(I/\phi J)\to\rMorph(\cC)
\qquad
(g:i\to\phi(j))\mapsto\mu_j\circ Fg
$$
(notation of \eqref{subsec_target-fctr} and \eqref{subsec_Morph-cat}).
We claim that $\tau^*$ factors through the map
$$
\Ob(I/\phi J)\to\Ob(I)
\quad :\quad
(g:i\to\phi(j))\mapsto i.
$$
Indeed, say that $g:i\to\phi(j)$ and $g':i\to\phi(j')$
are any two morphisms in $I$; we have to show that
$\tau^*_g=\tau^*_{g'}$. To this aim, since $i/\phi J$
is connected, we may assume that there exists
a morphism $h:g\to g'$ in $i/\phi J$; then
$$
\tau^*_g=\mu_j\circ Fg=\mu_{j'}\circ F\phi(h)\circ Fg=
\mu_{j'}\circ Fg'=\tau^*_{g'}
$$
as stated. Thus, let $\tau:\Ob(I)\to\rMorph(\cC)$
be the resulting map; we claim that $\tau$ is an element
of $L(X)$. Indeed, say that $h:i\to i'$ is any morphism
in $I$, and pick any $j,j'\in\Ob(J)$ with morphisms
$g:i\to\phi(j)$, $g':i'\to\phi(j')$ in $I$; we have
just seen that $\tau^*_{g'\circ h}=\tau^*_g$, which
translates as the identity : $\tau_i=\tau_{i'}\circ Fh$,
whence the claim. Lastly, by construction we have
$\beta_X(\tau)=\mu$, so $\beta_X$ is also surjective.
\end{proof}

\begin{remark}\label{rem_fun-cofinal}
(i)\ \
If $\phi:I\to J$ and $\psi:J\to K$ are cofinal functors,
the same holds for $\psi\circ\phi$. For the proof, we
may replace our universe $\sU$ by a larger one, after
which we may assume that $I,J,K$ are small categories,
and then the assertion follows directly from
proposition \ref{prop_MacL-cofinal}.

(ii)\ \
Let $\phi:J\to I$ be a functor between small categories.
Then $\phi$ is final if and only if it induces isomorphisms
of presheaves on $I$
\set\begin{equation}\label{eq_final}
\lim_IF\isom\lim_JF\circ\phi
\end{equation}
for every category $\cC$ with small $\Hom$-sets,
and every functor $F:I\to\cC$.

(iii)\ \
Let $\phi:J\to I$ and $F:I\to\cC$ be any two functors;
suppose that the induced morphism of presheaves
$\colim_\phi\one_\cC:\colim_IF\to\colim_JF\circ\phi$
is an isomorphism, and let $\tau:F\Rightarrow c_L$ be a
cocone. Then $\tau$ is universal if and only if the same
holds for $\tau*\phi:F\circ\phi\Rightarrow c_L$. Indeed,
$\tau$ and $\tau*\phi$ determine morphisms of presheaves
$t:h_{L^o}\to\colim_IF$ and $t':h_{L^o}\to\colim_JF\circ\phi$
on $\cC^o$ that assign to every $X\in\Ob(\cC)$ and every
section $f:L\to X$ of $h_{L^o}(X^o)$ the cocone
$\tau\odot c_f\in\colim_IF(X^o)$ and respectively the
cocone $(\tau*\phi)\odot c_f\in\colim_JF\circ\phi(X)$,
and $t$ (resp. $t'$) is an isomorphism if and only if
$\tau$ (resp. $\tau*\phi$) is universal.
However, $(\colim_\phi\one_\cC)\circ t$ is the morphism of
presheaves $h_{L^o}\isom\colim_JF\circ\phi$ that assigns to
every such $X$ and $f$ the cone
$(\tau\odot c_f)*\phi=(\tau*\phi)\odot c_f$, {\i.e.}
$t'=(\colim_\phi\one_\cC)\circ t$. So $t$ is an isomorphism
if and only if the same holds for $t'$, whence the contention.
Likewise, if \eqref{eq_final} is an isomorphism, a cone
$\tau':c_{L'}\Rightarrow F$ is universal if and only if the
same holds for $\tau'*\phi:c_{L'}\Rightarrow F\circ\phi$.
\end{remark}

\begin{lemma}\label{lem_filtered-final}
Let $J$ be a filtered category and $\phi:J\to I$ a functor.
We have :
\begin{enumerate}
\item
The functor $\phi$ is cofinal if and only if the following
two conditions hold :
\begin{enumerate}
\alphaenu
\item
For every $i\in\Ob(I)$ there exist $j\in\Ob(J)$ and
a morphism $i\to\phi(j)$ in $I$.
\item
For every $i\in\Ob(I)$, every $j\in\Ob(J)$ and every
pair of morphisms $f,g:i\to\phi(j)$ in $I$, there
exists a morphism $h:j\to j'$ in $J$ such that
$\phi(h)\circ f=\phi(h)\circ g$.
\end{enumerate}
\item
Moreover, if $\phi$ is cofinal, then $I$ is filtered.
\end{enumerate}
\end{lemma}
\begin{proof}(i): Indeed, suppose that the categories
$i/\phi J$ are connected for every $i\in\Ob(I)$; then
(a) clearly holds. To check that (b) holds as well, notice
that, since $J$ is filtered, $i/\phi J$ is locally directed
for every $i\in\Ob(I)$ (details left to the reader),
and therefore it is directed, by remark \ref{rem_cofinal}(i).
This implies that, for given $f,g$ as in (b), there
exist $j'\in\Ob(J)$ and morphisms $h_1,h_2:j\to j'$
such that $\phi(h_1)\circ f=\phi(h_2)\circ g$. But
since $J$ is filtered, we may then find a morphism
$h':j'\to j''$ in $J$ such that $h'\circ h_1=h'\circ h_2$,
so (b) holds with $h:=h'\circ h_1$.

Conversely, if (a) holds, then $i/\phi J$ is non-empty
for every $i\in\Ob(I)$. Next, let $g:i\to\phi(k)$ and
$g':i\to\phi(k')$ be any two objects of $i/\phi J$;
since $J$ is directed we may find morphisms $h:k\to j$
and $h':k'\to j$ for some $j\in\Ob(J)$, whence
objects $\phi(h)\circ g,\phi(h')\circ g':i\to\phi(j)$
of $i/\phi J$, and using (b) we deduce that there exist
an object $g'':i\to\phi(j')$ and morphisms $g\to g''$,
$g'\to g''$, {\em i.e.} $i/\phi J$ is directed.

(ii): Let us check that $I$ is directed : if $i$ and $i'$
are two objects of $I$, in view of condition (i.a), we may
find $j,j'\in\Ob(J)$ and morphisms $f:i\to\phi(j)$,
$f':i'\to\phi(j')$ in $I$, and since $J$ is directed, we
have as well morphisms $g:j\to j''$ and $g':j'\to j''$ in
$J$, for some $j''\in\Ob(J)$, whence morphisms
$\phi(g)\circ f:i\to\phi(k)$ and
$\phi(h')\circ f':i'\to\phi(k)$ in $I$. By remark
\ref{rem_cofinal}(i) it remains to check the coequalizing
condition of definition \ref{def_filtered-cols}(v), but the
latter is an immediate consequence of conditions (i.a) and (i.b)
\end{proof}

\begin{example}\label{ex_filtered-final}
(i)\ \ 
For instance, if $I$ is a filtered category, and $i$
is any  object of $I$, the category $i/I$ is again
filtered, and the target functor $\st_i:i/I\to I$ is
cofinal; indeed, one checks easily that $\st_i$
fulfills conditions (a) and (b) of lemma
\ref{lem_filtered-final}(i). Dually, if $I$ is cofiltered,
the category $I/i$ is again cofiltered, and the source
functor $\ss_i:I/i\to I$ is final.

(ii)\ \
If $I$ admits a final object $i_0$, then the (full)
subcategory $J$ of $I$ with $\Ob(J)=\{i_0\}$ fulfills
conditions (a) and (b) of lemma \ref{lem_filtered-final}(i),
so the inclusion functor $J\to I$ is cofinal (and $I$ is
trivially filtered). Then, let $F:I\to\cC$ be any functor,
and $\tau:F\Rightarrow c_L$ a universal cocone; taking into
account remark \ref{rem_fun-cofinal}(iii), we see that
$Fi_0$ represents the colimit of $F$, and $\tau_{i_0}:Fi\to L$
is an isomorphism in $\cC$. Dually, if $i'_0$ is any initial
object of $I$, and $\tau':c_{L'}\Rightarrow F$ any universal
cone, then $Fi'_0$ represents the limit of $F$ and
$\tau_{i'_0}:L'\to Fi'_0$ is an isomorphism.
\end{example}

\begin{example}\label{ex_fil-colim-in-Cat}
As an application, we give an explicit construction of
the filtered colimits of $\bCat$. Thus, let $I$ be a small
filtered category, and consider a functor
$$
\cC_\bullet:I\to\bCat
\qquad
i\mapsto\cC_i
\qquad
(\phi:i\to j)\mapsto(\cC_\phi:\cC_i\to\cC_j).
$$
We deduce a functor $\Ob(\cC_\bullet):I\to\Set$ that assigns
to every $i\in\Ob(I)$ the set $\Ob(\cC_i)$, and to every
morphism $\phi:i\to j$ in $I$ the map $\Ob(\cC_i)\to\Ob(\cC_j)$
defined by $\cC_\phi$. We set
$$
L:=\colim_I\Ob(\cC_\bullet).
$$
According to example \ref{ex_complete-cats}(i) -- and since
$I$ is filtered -- the elements of $L$ are the equivalence
classes $[i,X]$ of pairs $(i,X)$, where $i\in\Ob(I)$ and
$X\in\Ob(\cC_i)$, for the equivalence relation such that
$(i,X)\sim(j,Y)$ if and only if there exist $k\in\Ob(I)$
and morphisms $\phi_1:i\to k$ and $\phi_2:j\to k$ such that
$\cC_{\phi_1}X=\cC_{\phi_2}Y$. For every $i\in\Ob(I)$, let
$\st_i:i/I\to I$ be the target functor (see
\eqref{subsec_slice-cat}); notice that for every pair of
objects $i,j\in\Ob(I)$, the objects of the category
$(i,j)/I:=i/I\times_{(\st_i,\st_j)}j/I$ are
the pairs $i\xrightarrow{\phi_1}k\xleftarrow{\psi_2}j$ of
morphisms of $I$, and the morphisms
\set\begin{equation}\label{eq_morph-in-doublecomma}
(i,j)/\nu:(i\xrightarrow{\phi_1}k\xleftarrow{\phi_2}j)\to
(i\xrightarrow{\phi'_1}k'\xleftarrow{\phi'_2}j)
\end{equation}
are the morphisms $\nu:k\to k'$ of $I$ such that
$\nu\circ\phi_r=\phi'_r$ for $r=1,2$. For every couple
of pairs $((i,X),(j,Y))$ as in the foregoing,
we get a functor $h(i,X,j,Y):(i,j)/I\to\Set$ that assigns to
every $(i\xrightarrow{\phi_1}k\xleftarrow{\phi_2}j)\in\Ob((i,j)/I)$
the set $\Hom_{\cC_k}(\cC_{\phi_1}X,\cC_{\phi_2}Y)$, and to every morphism
$(i,j)/\nu$ as in \eqref{eq_morph-in-doublecomma} the induced map
$$
\Hom_{\cC_k}(\cC_{\phi_1}X,\cC_{\phi_2}Y)\to
\Hom_{\cC_{k'}}(\cC_{\phi'_1}X,\cC_{\phi'_2}Y)
\qquad
f\mapsto\cC_\nu(f)
$$
and we set
$$
H(i,X,j,Y):=\colim_{(i,j)/I}h(i,X,j,Y).
$$
Let also $\st_{i,j}:(i,j)/I\to I$ be the target functor given by
the rules : $(i\xrightarrow{\phi_1}k\xleftarrow{\phi_2}j)\mapsto k$,
and $(i,j)/\nu\mapsto\nu$; for every $i,j,k\in\Ob(I)$ we set as
well $(i,j,k)/I:=(i,j)/I\times_{(\st_{i,j},\st_k)}k/I$. Explicitly,
the objects of this category are all the systems
$\phi_\bullet:=(\phi_1:i\to t,\phi_2:j\to t,\phi_3:k\to t)$ of
morphisms of $I$ with a common target $t$ that we call {\em
the target of $\phi_\bullet$}. The morphisms
$$
(i,j,k)/\nu:\phi_\bullet\to\phi'_\bullet
$$
are the morphisms $\nu$ of $I$ with $\nu\circ\phi_r=\phi'_r$
for $r=1,2,3$. With this notation, for every third pair $(k,Z)$
as in the foregoing, we get a system of composition maps
$$
h(i,X,j,Y)(\phi_\bullet)\times h(j,Y,k,Z)(\phi_\bullet)\to
h(i,X,k,Z)(\phi_\bullet)
\qquad
\text{for every $\phi_\bullet\in\Ob((i,j,k)/I)$}
$$
given by the composition law of the category $\cC_t$,
where $t$ is the target of $\phi_\bullet$. Clearly this
system of maps is natural with respect to morphisms of
$(i,j,k)/I$; moreover we have obvious projection functors
$(i,k)/I\leftarrow(i,j,k)/I\to(i,j)/I$ and $(i,j,k)/I\to(j,k)/I$,
and using the criterion of lemma \ref{lem_filtered-final}(i)
it is easily seen all that these three functors are cofinal.
Therefore, taking colimits over $(i,j,k)/I$ of the foregoing
compatible system of composition maps yields a well defined map
$$
H(i,X,j,Y)\times H(j,Y,k,Z)\to H(i,X,k,Z)
\qquad
(f,g)\mapsto g\circ f
$$
and a simple inspection of the construction shows that
this composition law is associative : if $h\in H(k,Z,t,W)$
is any other element, we have
$h\circ(g\circ f)=(h\circ g)\circ f$ (details left to
the reader); likewise, the class of $\one_X$ in
$H(i,X,i,X)$ yields a left and right unit for this
composition law. Hence, let us choose for every $A\in L$
a representative $(i_A,X_A)$; we obtain a well defined
category $\cL$ with $\Ob(\cL)=L$ and with
$\Hom_\cL(A,B):=H(i_A,X_A,i_B,X_B)$ for every $A,B\in L$,
with the composition law thus defined. Moreover, for every
$i\in I$ we obtain a well defined functor $G_i:\cC_i\to\cL$
by setting $GX:=[i,X]$ for every $X\in\Ob(\cC_i)$, and
by assigning to every morphism $X\to Y$ of $\cC_i$ its
class in $H(i,X,i,Y)$. Lastly, let
$F:\cC_\bullet\Rightarrow\cD$ be a cocone with basis
$\cC_\bullet$ and whose vertex is any category $\cD$.
The induced cocone $\Ob(\cC_\bullet)\Rightarrow\Ob(\cD)$
yields a well defined map $L\to\Ob(\cD)$. Moreover,
for every two pairs $(i,X)$ and $(j,Y)$ and every object
$(i\xrightarrow{\phi_1}k\xleftarrow{\phi_2}j)$ of
$(i,j)/I$, the functor $F_k$ yields a map
$$
\Hom_{\cC_k}(\cC_{\phi_1}X,\cC_{\phi_2}Y)\to
\Hom_\cD(F_k\circ\cC_{\phi_1}X,F_k\circ\cC_{\phi_2}Y)=
\Hom_\cD(F_iX,F_jY)
$$
which defines a cocone with vertex $\Hom_\cD(F_iX,F_jY)$
and for basis the functor $h(i,X,j,Y)$, whence a
well defined map
$$
F_{i,X,j,Y}:H(i,X,j,Y)\to\Hom_\cD(F_iX,F_jY).
$$
A simple inspection shows that
$F_{i,X,k,Z}(g\circ f)=F_{j,Y,k,Z}(g)\circ F_{i,X,j,Y}(f)$ for
every $f,g$ as in the foregoing (details left to the reader).
We obtain therefore a well defined functor $F:\cL\to\cD$ by
the rule : $A\mapsto F_{i_A}X_A$ and $F_{AB}f:=G_{i_A,X_A,i_B,X_B}$
for every $A,B\in L$ and every $f\in\Hom_\cL(A,B)$. By
construction, we have $F\circ G_i=F_i$ for every $i\in I$,
and this is clearly the unique functor fulfilling this
system of identities, so the proof is concluded.
\end{example}

\begin{remark}\label{rem_cofinally-small}
(i)\ \
Let $I$ be a finally small category, $\cC$ a category
with small $\Hom$-sets, and $F:I\to\cC$ any functor.
By definition, there exists a final functor $\phi:J\to I$
with $J$ a small category, and therefore we may define
the presheaf
$$
L_{F,\phi}:=\lim_JF\circ\phi.
$$
We claim that $L_{F,\phi}$ is independent - up to natural
isomorphism - of the choice of $\phi$. To see this, let
us choose a universe $\sU'$ such that $\sU\subset\sU'$
and such that $I$ is a small $\sU'$-category. By
remark \ref{rem_fun-cofinal}(i), the induced morphism
of $\sU'$-presheaves
$$
\omega_\phi:\lim_IF\to L_{F,\phi}
$$
is an isomorphism. Therefore, if $\phi':J'\to I$ is
any other final functor with $J'$ also small, we get
an isomorphism of $\sU'$-presheaves
$\omega_{\phi'}\circ\omega_\phi^{-1}:L_{F,\phi}\isom L_{F,\phi'}$.
Since the inclusion \eqref{eq_change-universe} is
fully faithful, the latter is actually an isomorphism
of $\sU$-presheaves, as required.

(ii)\ \
In view of (i), we may define the {\em limit of $F$}
as $L_{F,\phi}$, for any choice of $\phi$. This presheaf
is therefore well defined, up to natural isomorphism, and
we denote it as usual by
$$
\lim_IF.
$$
We also say that $\lim_IF$ is a {\em finally small limit}
of $\cC$. Of course, if $I$ is already small, we can choose
$\phi:=\one_I$, so the above notation is compatible with
that of definition \ref{def_limits}(i).

(iii)\ \
Dually, if $I$ is cofinally small, we can define
the {\em colimit of $F$}
$$
\colim_IF:=\colim_JF\circ\phi
$$
for any choice of cofinal functor $\phi:J\to I$ with
$J$ small, and again the resulting presheaf is well
defined up to natural isomorphism, and we also say
that it is a {\em cofinally small colimit} of $\cC$.
Clearly, if $\cC$ is complete (resp. cocomplete), then
every finally small (resp. cofinally small) limit is
representable in $\cC$. As an application, we may
prove Freyd's {\em adjoint functor theorem}, which is
the following partial converse of proposition
\ref{prop_was-get-maddd}(iii) :
\end{remark}

\begin{theorem}\label{th_adjoint-fctr-th}
Let $\cA$ be a complete category, $F:\cA\to\cB$ a functor,
and suppose that $\cA$ and $\cB$ have small $\Hom$-sets
(see \eqref{subsec_opposing}). The following conditions
are equivalent :
\begin{enumerate}
\alphaenu
\item
$F$ admits a left adjoint.
\item
$F$ commutes with all the limits of $\cA$, and every object
$B$ of $\cB$ admits a {\em solution set}, {\em i.e.} an essentially
small subset $S_B\subset\Ob(\cA)$ such that, for every $A\in\Ob(\cA)$,
every morphism $f:B\to FA$ admits a factorization of the form
$$
f=Fh\circ g
$$
for $h:A'\to A$ a morphism in $\cA$ with $A'\in S_B$,
and $g:B\to FA'$ a morphism in $\cB$.
\end{enumerate}
\end{theorem}
\begin{proof} If $F$ admits a left adjoint $G$, then $F$
commutes with all limits of $\cA$, by proposition
\ref{prop_was-get-maddd}(iii), and it is easily seen that
$\{GB\}$ is a solution set for every object $B\in\Ob(\cB)$.

Conversely, fix any $B\in\Ob(\cB)$; under the stated
assumptions, example \ref{ex_cofiltered-comma}(i) says
that the category $B/F\cA$ is cofiltered. Denote by
$\cE_B$ the full subcategory of $B/F\cA$ whose objects
are all the morphisms $B\to FA$ with $A\in S_B$. Then
condition (b) means that for every object $X$ of
$B/F\cA$ there exists an object $X'\in\cE_B$ with
a morphism $X'\to X$ of $B/F\cA$. It follows easily
that $\cE_B$ is cofiltered as well; then the
inclusion functor $\cE_B\to B/F\cA$ is final, by
lemma \ref{lem_filtered-final}(i), and $\cE_B$ is
small, hence $B/F\cA$ is finally small. Following
remark \ref{rem_cofinally-small}(ii), we may thus
define
$$
G'B:=\lim_{B/F\cA}\ss_B
\qquad
G'f:=\lim_{f/F\cA}\one_\cA
$$
for any $B\in\Ob(\cB)$ and any morphism $f:B'\to B$
(notation of \eqref{subsec_fibreovercat}). Next, since
$\cA$ is complete, $G'$ is representable by some
functor $G:\cB\to\cA$ (remark \ref{rem_represent-morph}(ii)),
and we claim that $G$ is the sought left adjoint. To check
this assertion, we may replace the universe $\sU$ by a
larger one, after which we may assume that $\cA$ and $\cB$
are small, and therefore $G'B$ is the usual limit of
definition \ref{def_limits}(i). We shall then exhibit
explicit unit and counit for the adjoint pair $(G,F)$,
as follows. Let $\tau^B:c_{GB}\Rightarrow\ss_B$ be a
universal cone; since $F$ commutes with limits, the cone
$F*\tau^B:c_{FGB}\Rightarrow F\circ\ss_B$ is still
universal (remark \ref{rem_representable-case}(i));
on the other hand, the functor $B/F:B/F\cA\to B/\cB$
can be regarded as a cone $B/F:c_B\Rightarrow F\circ\ss_B$,
so there exists a unique morphism $\eta_B:B\to FGB$ such
that $B/F=(F*\tau^B)\odot c_{\eta_B}$. Likewise, for
any $A\in\Ob(\cA)$ we get a morphism
$\eps_A:=\tau^{FA}_{(FA,\one_{FA})}:GFA\to A$. It is easily
seen that $\eta$ (resp. $\eps$) is a natural transformation
$\one_\cB\Rightarrow FG$ (resp. $GF\Rightarrow\one_\cA$),
and the identity $(F*\eps)\odot(\eta*F)=\one_F$ is immediate
from the construction. Lastly, we fix $B\in\Ob(\cB)$ and we
check that $(\eps*G)_B\circ(G*\eta)_B=\one_{G(B)}$; to this aim,
and due the universality of $\tau^B$, it suffices to show that
$$
\tau_\psi^B=\tau^B_\psi\circ\eps_{GB}\circ G(\eta_B)
\qquad
\text{for every object $\psi:B\to FA$ of $B/F\cA$}.
$$
However, notice that the construction of $Gf$ yields
the identity
\set\begin{equation}\label{eq_tricky-Freyd}
\tau^{B'}_\psi\circ Gf=\tau^B_{\psi\circ f}
\qquad
\text{for all morphisms $f:B\to B'$ and $\psi:B'\to FA$ in $\cC$}.
\end{equation}
Thus, we may compute :
$$
\tau^B_\psi\circ\eps_{GB}\circ G(\eta_B)
=\tau^B_\psi\circ\tau^{FGB}_{\one_{FGB}}\circ G(\eta_B)
=\tau^{FGB}_{F(\tau^B_\psi)}\circ G(\eta_B)
=\tau^B_{F(\tau^B_\psi)\circ\eta_B}
= \tau^B_\psi
$$
where the second equality holds by the naturality of
$\tau^{FGB}$, the third follows from \eqref{eq_tricky-Freyd},
and the fourth comes from the construction of $\eta_B$.
\end{proof}

Of course, the ``dual'' of theorem \ref{th_adjoint-fctr-th}
yields a criterion for the existence of right adjoints.

\begin{example}\label{ex_lim_interchange}
(i)\ \
Let $J$ and $K$ be any two small categories; set
$I:=J\times K$ and let $\pi:I\to J$ be the projection
functor. For every $j\in\Ob(J)$ we have an induced functor
$$
\iota_j:K\to j/\pi I
\qquad
k\mapsto((j,k),\one_j)
\qquad
(g:k\to k')\mapsto(\one_j,g)
$$
and we claim that $\iota_j$ is final (definition
\ref{def_MacL-cofinal}(ii)). Indeed, for any
$((j_0,k_0),f_0:j\to j_0)\in\Ob(j/\pi I)$, the
category $\iota_j K/((j_0,k_0),f_0)$ is isomorphic to
$K/k_0$, which is obviously connected, whence the claim.

(ii)\ \
Moreover, every morphism $f:j\to j'$ in $J$ induces
a natural transformation
$$
\beta^f:\iota_j\Rightarrow(f/\pi I)\circ\iota_{j'}
$$
which assigns to every $k\in\Ob(K)$ the morphism
$(f,\one_k):((j,k),\one_j)\to((j',k),f)$. Then,
$\st_j*\beta^f$ is a natural transformation
$\st_j\circ\iota_j\Rightarrow
\st_j\circ(f/\pi I)\circ\iota_{j'}=\st_{j'}\circ\iota_{j'}$
and a direct inspection of the definitions yields a
commutative diagram for every functor $F:I\to\cC$
$$
\xymatrix{
\int^\wedge_\pi F(j) \ar[rrr]^-{\int^\wedge_\pi F(f)} \ar[d] & & &
\int^\wedge_\pi F(j') \ar[d] \\
\lim_K F\circ\st_j\circ\iota_j \ar[rrr]^-{\lim_K(F*\st_j*\beta^f)}
& & & \lim_K F\circ\st_{j'}\circ\iota_{j'}
}$$
whose vertical arrows are the natural isomorphisms
provided by (i) and remark \ref{rem_fun-cofinal}(ii) and
whose bottom arrow is given by remark \ref{rem_wishful}(ii).
In other words, the functor $\int^\wedge_\pi$ is naturally
isomorphic to the one that assigns to every such $F$
the functor
$$
\int^\wedge_KF:J\to\cC^\wedge
\qquad
j\mapsto\lim_K F\circ\st_j\circ\iota_j
\qquad
(f:j\to j')\mapsto\lim_K(F*\st_j*\beta^f).
$$
Lastly, proposition \ref{prop_Fubini} can be restated
in this case more synthetically, as the existence of
a natural isomorphism in $\cC^\wedge$
$$
\lim_{J\times K}F\isom\Lim_J\int^\wedge_KF
\qquad
\text{for every functor $F:J\times K\to\cC$}.
$$

(iii)\ \
On the other hand, we may also apply the same considerations
to the second projection functor $I\to K$. Summing up, we
deduce a natural isomorphism
$$
\Lim_J\int^\wedge_K F\isom\Lim_K\int^\wedge_J F
\qquad
\text{for every functor $F:J\times K\to\cC$}
$$
which expresses the well known {\em interchange property}
for double limits.

(iv)\ \
Dually, we obtain as well an interchange property for double
colimits, that states the existence of natural isomorphisms
$$
\Colim_J\int_\wedge^K F\isom\Colim_{J\times K}F\isom
\Colim_K\int_\wedge^J F
\qquad
\text{for every functor $F:J\times K\to\cC$}
$$
where $\int_\wedge^K$ is the opposite of the functor
$\int^\wedge_{K^o}$ (and likewise for $\int^J_\wedge$ :
details left to the reader).
\end{example}

\begin{remark}\label{rem_interchange}
(i)\ \
Keep the situation of example \ref{ex_lim_interchange},
and suppose moreover that $\cC$ is complete (resp.
cocomplete). Then the functor $\int^\wedge_K$ (resp.
$\int^K_\wedge$) is isomorphic to the composition of
$\bFun(J,h_\cC)$ (resp. of $\bFun(J,h^o_{\cC^o})$)
and a functor
$$
\int_K:\bFun(I,\cC)\to\bFun(J,\cC)
\qquad
\text{(resp. $\displaystyle\int^K:\bFun(I,\cC)\to\bFun(J,\cC)$)}
$$
that is well defined up to isomorphism, and is isomorphic
to $\int_\pi$ (resp. to $\int^\pi$), so the interchange
properties of example \ref{ex_lim_interchange}(iii,iv)
can, in this case, be expressed also in terms of these
new functors (details left to the reader).

(ii)\ \
Furthermore, if $\cC$ is both complete and cocomplete,
we get a natural morphism
\set\begin{equation}\label{eq_interchange}
\Colim_J\int_K F\to\Lim_K\int^J F
\qquad
\text{for every functor $F:J\times K\to\cC$}
\end{equation}
as follows. First, for every $(j,k)\in\Ob(J\times K)$ we
have a natural morphism
$$
\omega_{j,k}:\int_K F(j)\to F(j,k)\to\int^J F(k)
$$
given by the choice of a universal cone and cocone for
the limit and colimit that are represented respectively
by the source and target of this morphism. By inspecting
the constructions, we then check easily that, for fixed
$j\in\Ob(J)$, the system $(\omega_{j,k}~|~k\in\Ob(K))$
is a cone with vertex $\int_KF(j)$; the choice of a
universal cone for the functor $\int^KF$ then determines
a unique morphism
$$
\omega_j:\int_KF(j)\to\Lim_J\int^KF
\qquad
\text{for every $j\in\Ob(J)$}.
$$
In turns, the system $(\omega_j~|~j\in\Ob(J))$ is a
cocone which determines \eqref{eq_interchange}, after
fixing a universal cocone for the functor $\int_KF$.
Of course, \eqref{eq_interchange} depends on the choices
of all these universal cones and cocones, but two different
sets of such choices will modify the morphism by left
and right composition with uniquely determined isomorphisms.

(iii)\ \
The question rises, whether \eqref{eq_interchange} is
an isomorphism. This is not always the case. If
\eqref{eq_interchange} is an isomorphism for every
functor $F:J\times K\to\cC$, we say that in the
category $\cC$ {\em the limits indexed by $K$ commute
with the colimits indexed by $J$}.

(iv)\ \
For instance, if $K$ is finite and $J$ is filtered,
then the limits indexed by $K$ commute with the colimits
indexed by $J$ in the category $\Set$ (\cite[Th.2.13.4]{Bor}).
We say briefly that {\em the finite limits commute with the
filtered colimits} in $\Set$.

Taking into account example \ref{ex_big-complete}(i),
we deduce more generally that, for any category
$\cC$, the finite limits in $\cC^\wedge$ commute with all
filtered colimits.
\end{remark}

Here is another useful application to Kan extensions :

\begin{proposition}\label{prop_faithful-Kan-ext}
In the situation of theorem {\em\ref{th_Kan-ext}}, suppose
that $\phi$ is fully faithful. Then the same holds for the
right (resp. left) Kan extension along $\phi$.
\end{proposition}
\begin{proof} Consider the right Kan extension along $\phi$
(and hence, we assume that the relevant completeness conditions
for $\cC$ as in remark \ref{rem_Lims-and-Colims}(iii,iv) are
fulfilled). Let $F:I\to\cC$ be any functor, and to ease
notation set $K:=\int_\phi F$; by proposition
\ref{prop_fullfaith-adjts}(iii), it suffices to exhibit
an isomorphism of functors $\eps:K\circ\phi\isom F$.
Now, recall that $Kj\in\cC$ represents the limit of the functor
$F\circ\st_j:j/\phi I\to\cC$, for every $j\in\Ob(J)$, and
the construction of $K$ involves the choice of a universal
cone $\tau^j:c_{Kj}\Rightarrow F\circ\st_j$ for every such $j$.
However, notice that $(i,\one_{\phi(i)})$ is an initial object
of $\phi(i)/\phi I$ : indeed, let $(i',f:\phi(i)\to\phi(i'))$
be any other object; since $\phi$ is fully faithful, there
exists a unique morphism $g:i\to i'$ in $J$ such that
$\phi(g)=f$, and then $(g,\one_f)$ is the unique morphism
$(i,\one_{\phi(i)})\to(i',f)$ in $\phi(i)/\phi I$. By example
\ref{ex_filtered-final}(ii), it follows that
$$
\eps_i:=\tau^{\phi(i)}_{(i,\one_{\phi(i)})}:K\circ\phi(i)\to Fi
$$
is an isomorphism for every $i\in\Ob(I)$. It remains to check
the naturality of the rule : $i\mapsto\eps_i$. Hence, let
$g:i\to i'$ be any morphism of $I$; by construction, for every
$X\in\Ob(\cC)$ we have a commutative diagram of sets :
$$
\xymatrix{ h_{K\phi(i)}(X) \ar[r] \ar[d]_{h_{K\phi(g)}} &
\lim_{\phi(i)/\phi I}F\circ\st_{\phi(i)}(X) \ar[d] \\
h_{K\phi(i')}(X) \ar[r] &
\lim_{\phi(i')/\phi I}F\circ\st_{\phi(i')}(X)
}$$
where $h$ denotes as usual the Yoneda embedding, and where
the top and bottom horizontal arrows are the bijections
induced by $\tau^i$ and $\tau^{i'}$ respectively. The right
vertical arrow is the map that sends every cone
$\beta:c_X\Rightarrow F\circ\st_{\phi(i)}$ to the cone
$\beta*(\phi g/\phi I):c_X\Rightarrow F\circ\st_{\phi(i')}$,
where $\phi g/\phi I:\phi(i')/\phi I\to\phi(i)/\phi I$ is
the functor defined as in \eqref{subsec_fibreovercat}.
Taking $X:=K\phi(i)$, we see that the top horizontal (resp.
left vertical) arrow maps $\one_{K\phi(i)}$ to $\tau^{\phi(i)}$
(resp. to $K\phi(g)$); then we get the identity :
$$
\tau^{\phi(i')}\odot c_{K\phi(g)}=\tau^{\phi(i)}*(\phi g/\phi I)
$$
whence
$\eps^{i'}\circ K\phi(g)=\tau^{\phi(i)}_{(i',\phi(g))}:K\phi(i)\to Fi'$.
Lastly, the morphism $g$ induces a morphism
$\phi(i)/g:(i,\one_{\phi(i)})\to(i',\phi(g))$ of
$\phi(i)/\phi I$, whence the identity :
$\tau^{\phi(i)}_{(i',\phi(g))}=Fg\circ\tau^{\phi(i)}_{(i,\one_{\phi(i)})}$.
Summing up, we conclude that
$$
\eps^{i'}\circ K\phi(g)=Fg\circ\eps^i
$$
as required. The assertion for left Kan extensions admits
the dual proof.
\end{proof}

\begin{corollary}\label{cor_lable}
Let $\cB$, $\cC$ be two small categories, and $F:\cB\to\cC$
a functor. We have :
\begin{enumerate}
\item
If $\cB$ is finitely complete and $F$ is left exact, 
then $F_!$ is exact.
\item
$F$ is fully faithful $\Leftrightarrow$ $F_!$ is fully faithful
$\Leftrightarrow$ $F_*$ is fully faithful.
\end{enumerate}
\end{corollary}
\begin{proof}(i): Under these assumptions, the category
$Y/F\cB$ is cofiltered for every $Y\in\Ob(\cC)$ (example
\ref{ex_cofiltered-comma}(i)), hence $(Y/F\cB)^o$ is
filtered. However, the filtered colimits in the category
$\Set$ commute with all finite limits (see remark
\ref{rem_interchange}(iv)), so the assertion follows
from remark \ref{rem_was-cofinal}(iii) and proposition
\ref{prop_was-get-maddd}(iv).

(ii): By proposition \ref{prop_faithful-Kan-ext} we know
that if $F$ is fully faithful, the same holds for $F_!$
and $F_*$. Also, by proposition \ref{prop_fullfaith-adjts}(iv),
$F_!$ is fully faithful if and only if the same holds for $F_*$.
Lastly,  if $F_!$ is fully faithful, then \eqref{eq_ess-comm-dig}
and the full faithfulness of the Yoneda embeddings, imply that
$F$ is fully faithful.
\end{proof}

\begin{example}\label{ex_Cat-is-coco}
Let $I$ be any small category; we wish to apply Freyd's
adjoint theorem \ref{th_adjoint-fctr-th} in order to
exhibit a left adjoint
$$
\Colim_I:\bFun(I,\bCat)\to\bCat
$$
for the corresponding constant functor $c$ for the category
$\bCat$. In particular, this will prove that all the
colimits indexed by $I$ are representable in $\bCat$, and
since $I$ is arbitrary, we will conclude that $\bCat$
is cocomplete. Now, we know already that $\bCat$ is
complete (example \ref{ex_cat-cats}(i)), and $c$ commutes
with all limits of $\bCat$, by virtue of corollary
\ref{lem_presheaf-in-a-cat}(ii). It remains to check
that every functor $\cC_\bullet:I\to\bCat$ admits a
solution set $S\subset\Ob(\bCat)$. To this aim, let
$\cD:=:=\amalg_{i\in\Ob(I)}\cC_i$ be the coproduct of the
family of categories $(\cC_i~|~i\in\Ob(I))$, as in
example \ref{ex_cat-cats}(i); for every small category
$\cB$ and every co-cone $F_\bullet:\cC_\bullet\Rightarrow\cB$,
let $F:\cD\to\cB$ be the resulting functor, and denote
by $\Img(F)\subset\cB$ the image of $F$, defined as in
example \ref{ex_pre-misc} (it is easily seen that
$\bCat$ is well-powered, so the image is well defined).
The system of categories
$$
\cF:=(\Img(F)~|~\cB\in\Ob(\bCat),F_\bullet:\cC_\bullet\to\cB)
$$
is not small, but if we pick in $\cF$ a representative
for each isomorphism class, we do get a small family
(details left to the reader); it is easily seen that
any such choice of representatives yields a solution
set for $\cC_\bullet$.
\end{example}

In applications, often the indexing category of a
limit (or colimit) is a partially ordered set, and
usually such limits are easier to handle than limits
over general indexing categories. Now, given such a
general indexing category $I$, one may try to find
a cofinal functor $J\to I$ from a partially ordered
set $J$. The following proposition says that this
can always be achieved, at least if $I$ is filtered.

\begin{proposition}\label{prop_filter-Deligne}
Let $I$ be any filtered category. The following holds :
\begin{enumerate}
\item
There exists a cofinal functor $J\to I$, with $J$ a
filtered partially ordered set.
\item
If\/ $I$ is small, we can find $J$ as in {\em(i)},
such that $J$ is also small.
\end{enumerate}
\end{proposition}
\begin{proof} Let us say that $I$ admits a
{\em largest object}, if there exists $i\in\Ob(I)$
such that $\Hom_I(i',i)\neq\emptyset$ for every
$i'\in\Ob(I)$. We notice :

\begin{claim}\label{cl_Deligne-0}
Let $j\in\Ob(I)$ be any object, and
$(f_t:i_t\to j~|~t=1,\dots,n)$ any finite system of
morphisms in $I$. If $I$ does not have a largest
object, there exists a morphism $j\to i$ in $I$ such
that $i\neq i_t$ for every $t=1,\dots,n$.
\end{claim}
\begin{pfclaim} Suppose the claim fails, in which
case we may assume that $j=i_1$. Since $I$
is filtered, for every $k\in\Ob(I)$ we may find
morphisms $g:k\to k'$ and $j\to k'$, and the
assumption implies that $k'=i_{n(k)}$ for some
$n(k)\in\{1,\dots,n\}$. Set $h:=f_{n(k)}\circ g$,
so $h$ is a morphism $k\to i_1$, and especially
$\Hom_I(k,i_1)\neq\emptyset$ for every $k\in\Ob(I)$,
which is absurd, since $I$ does not have a largest
object.
\end{pfclaim}

Next, let $\N$ be the filtered category associated
with the totally ordered set of natural numbers; we
remark that if $I$ is any filtered category, then
$\N\times I$ is still filtered, and it does not
admit a largest object; moreover, one checks easily
that the projection functor $\N\times I\to I$
fulfills conditions (a) and (b) of lemma
\ref{lem_filtered-final}(i), so it is cofinal.
In light of remark \ref{rem_fun-cofinal}(i), we may
then replace $I$ by $\N\times I$, and assume from
start that $I$ does not have a largest object.

Let us say that a {\em diagram} of $I$ is any pair
$D:=(A,B)$ with $A\subset\Ob(I)$,
$B\subset\rMorph(I)$, and such that, for
every $f\in B$, the source and target of $f$ lie
in $A$. We say that $A$ (resp. $B$) is the
set of {\em objects} (resp. {\em morphisms}) of $D$.
An element $e$ of $A$ is said to be {\em final}
in $(A,B)$ if the following holds :
\begin{itemize}
\item
for every $i\in\Ob(I)$, the set $\Hom_I(i,e)\cap B$
contains exactly one element $f^D_i$
\item
$f^D_e=\one_e$, and for every morphism $g:i\to j$ in $B$
we have $f^D_j\circ g=f^D_i$.
\end{itemize}
Denote by $J$ the set consisting of all diagrams
$(A,B)$ which admit a unique final element $e(A,B)$,
and such that $A$ and $B$ are finite sets. We order
$J$ by inclusion, {\em i.e.} $(A,B)\leq(A',B')$ if
and only if $A\subset A'$ and $B\subset B'$. Then,
if $D,D'\in J$ and $D\leq D'$, there exists a
unique morphism $e(D,D'):e(D)\to e(D')$ in the
set of morphisms of $D'$, and if $D\leq D'\leq D''$,
then clearly $e(D,D'')=e(D',D'')\circ e(D,D')$.
Thus, the rule $D\mapsto e(D)$ yields a well
defined functor $e:J\to I$, and we have to prove
that $e$ is cofinal.

To this aim, we check that conditions (a) and
(b) of lemma \ref{lem_filtered-final}(i) hold
for $e$. Now, if $i$ is any object of $I$, the
diagram $D_i:=(\{i\},\{\one_i\})$ lies in $J$,
and $\one_i:i\to e(D_i)$ is a morphism in $I$,
so condition (a) holds. Next, let $D=(A,B)\in J$
be any diagram, $i\in\Ob(I)$ any object, and
$f,g:i\to e(D)$ any two morphisms in $I$;
we notice :

\begin{claim}\label{cl_Deligne}
There exist $i'\in\Ob(I)\setminus A$ and a morphism
$h:e(D)\to i'$ such that $h\circ f=h\circ g$.
\end{claim}
\begin{pfclaim} Since $I$ is filtered, we can
find a morphism $h:e(D)\to i'$ that coequalizes $f$
and $g$, and claim \ref{cl_Deligne-0} says that we may
choose $i\notin A$.
\end{pfclaim}

Now, let $i'$ be as in claim \ref{cl_Deligne},
and $D'$ the diagram of $I$ whose set of objects
is $A\cup\{i'\}$, and whose set of morphisms is
$B\cup\{h\circ f_j~|~j\in A\}$. It is easily seen
that $D'\in J$, and by construction, $e(D')=i'$
and $e(D,D')=h$, so (b) holds.

It remains to show that $J$ is filtered. However,
say that $D=(A,B)$ and $D'=(A',B')$ are any two
elements of $J$; we notice :

\begin{claim}\label{cl_Deligne-2}
There exist $i\in\Ob(I)\!\setminus\!(A\cup A')$
and morphisms $f:e(D)\to i$ and $f':e(D')\to i$
in $I$.
\end{claim}
\begin{pfclaim} This is similar to the proof of
claim \ref{cl_Deligne} : we can always find
$f:e(D)\to i$ and $f':e(D')\to i$, and by claim
\ref{cl_Deligne-0} we may arrange that
$i\notin A\cup A'$.
\end{pfclaim}

Hence, let $f,f'$ be as in claim \ref{cl_Deligne-2},
and denote $D(f,f')$ the diagram whose set of
objects is $A\cup A'\cup\{i\}$, and whose set of
morphisms is $B\cup B'\cup\{f\circ f^D_j~|~j\in A\}\cup
\{f'\circ f^{D'}_{j'}~|~j'\in A'\}$. Clearly, the
assertion will follow from :

\begin{claim} We can choose $f$ and $f'$ so that
$D(f,f')\in J$.
\end{claim}
\begin{pfclaim} It amounts to checking that,
for suitable $f$ and $f'$, we have
$f\circ f^D_j=f'\circ f^D_j$ for every $j\in A\cap A'$.
However, if $f$ and $f'$ as in claim \ref{cl_Deligne-2},
have been picked, a simple induction on the cardinality
of $A\cap A'$ shows that we may find $h:i\to i'$ such
that
$$
h\circ f\circ f^D_j=h\circ f'\circ f^D_j
\qquad
\text{for every $j\in A\cap A'$}
$$
and by claim \ref{cl_Deligne-0}, we may arrange
that $i'\notin A\cup A'$. Then the diagram
$D(h\circ f,h\circ f')$ will do.
\end{pfclaim}

Lastly, we remark that if $I$ is small, the same holds
for $J$.
\end{proof}

\begin{example}\label{ex_count-filtered}
Let us say that a category $I$ is {\em countable} if both
$\Ob(I)$ and $\rMorph(I)$ are countable sets. A simple
inspection of the proof of proposition
\ref{prop_filter-Deligne} shows that if $I$ is a filtered
countable category, then there exist a countable filtered
partially ordered set $J:=\{j_n~|~n\in\N\}$ and a cofinal
functor $J\to I$. But it is easy to construct a cofinal
functor $f:\N\to J$ (where $\N$ is endowed with its standard
total ordering) : indeed, we may choose $f(0):=j_0$, and then
pick inductively for every $n\in\N$ an element $f(n+1)\in J$
larger than $j_{n+1}$ and $f(n)$. Summing up, we see
that for every countable filtered category $I$ there exists
a cofinal functor $\N\to I$.
\end{example}

\subsection{Localizations of categories}
\label{subsec_graphs}
Given a category $\cC$ and an arbitrary set of morphisms
$\Sigma\subset\rMorph(\cC)$, we wish to describe
a general procedure that adds formally to $\cC$ an inverse
$f^{-1}$ for every $f\in\Sigma$. More details can be found
in \cite[Ch.5]{Bor}. We begin with the following :

\begin{definition}\label{def_graph}
(i)\ \ 
A {\em graph} $\cG$ is the datum of 
\begin{itemize}
\item
a set $V(\cG)$, whose elements are called the {\em vertices}
of $\cG$
\item
for every $A,B\in V(\cG)$ a set $\cG(A,B)$, whose elements
are called the {\em arrows} from $A$ to $B$. If $f\in\cG(A,B)$,
we say also that $A$ is the {\em source} and $B$ is the
{\em target} of $f$.
\end{itemize}

(ii)\ \
A {\em morphism} of graphs $F:\cG\to\cG'$ consists of
\begin{itemize}
\item
A map of sets $F:V(\cG)\to V(\cG')$
\item
For every $A,B\in V(\cG)$, a map of sets
$\cG(A,B)\to\cG'(FA,FB)$ : $f\mapsto Ff$.
\end{itemize}

(iii)\ \
A graph is {\em small} if $V(\cG)$ is a small set, and
$\cG(A,B)$ is a small set for all $A,B\in V(\cG)$.
\end{definition}

Obviously, morphisms of graphs can be composed. so we
have a category
$$
\sU\tdu\bGraph
$$
whose objects are all the small graphs. As usual, we shall
omit the prefix $\sU$, unless we have to deal with more
than one universe. Clearly, there is a natural functor
\set\begin{equation}\label{eq_cat-to-graph}
\bCat\to\bGraph
\end{equation}
which assigns to any category $\cC$ its underlying graph,
whose vertices are the objects of $\cC$, and whose arrows
are the morphisms of $\cC$, and simply forgets the composition
law of $\cC$. We are going to construct an explicit
left adjoint for the functor \eqref{eq_cat-to-graph}.
To this aim, let $\cG$ be any graph and $n\in\N$ any
integer; a {\em path of length $n$} in $\cG$ is a
sequence
$$
p:=(A_0,f_1,A_1,f_2,\dots,A_n)
$$
alternating vertices $A_0,\dots,A_n$ of $\cG$ and
arrows $f_1,\dots,f_n$, such that the source and
target of $f_i$ are respectively $A_{i-1}$ and $A_i$
for every $i=1,\dots,n$. Then $A_0$ is called the
{\em source} and $A_n$ the {\em target} of $p$. Given
paths $p:=(A_0,f_1,A_1,\dots,A_n)$ and
$p':=(A_n,f_{n+1},A_{n+1},\dots,A_{n+m})$ of lengths
respectively $n$ and $m$, and such that the target
of $p$ equals the source of $p'$, we obtain a path
of length $n+m$, by the rule :
$$
p'\circ p:=(A_0,f_1,\dots,A_n,f_n,\dots,A_{n+m}).
$$
Denote by $P(\cG)$ the set of all paths of $\cG$
(of arbitrary length); with this composition law,
clearly $P(\cG)$ is the set of morphisms of a
{\em path category} whose set of objects is $V(\cG)$
(for any $A\in V(\cG)$, the path $(A)$ of length zero
serves as the identity endomorphism of $A$). Also,
if $F:\cG\to\cG'$ is any morphism of graphs, we have
an induced functor
$$
P(F):(V(\cG),P(\cG))\to(V(\cG'),P(\cG'))
$$
whose map on objects is $F$, and such that
$P(F)p:=(FA_0,Ff_1,\dots,FA_n)$ for every path
$p:=(A_0,f_1,\dots,A_n)$ of $\cG$. We have thus
defined a functor
\set\begin{equation}\label{eq_graph-to-cat}
\bGraph\to\bCat
\qquad
\cG\mapsto(V(\cG),P(\cG)).
\end{equation}

\begin{proposition}\label{prop_graph-to-cat}
The functor \eqref{eq_graph-to-cat} is left adjoint
to the forgetful functor
\eqref{eq_cat-to-graph}.
\end{proposition}
\begin{proof} Indeed, given a graph $\cG$, a category
$\cC$ and a morphism of graphs $F:\cG\to\cC$, define
a functor $G:(V(\cG),P(\cG))\to\cC$ by the rule
$$
GA:=FA
\qquad
Gp:=Ff_n\circ\cdots Ff_1
$$
for every $A\in V(\cG)$ and every
$p:=(A_0,f_1,\dots,A_n)\in P(\cG)$. Clearly $G$ is the
unique functor whose restriction to $\cG$ agrees with
$F$, whence the contention.
\end{proof}

\begin{corollary}\label{cor_char-finally-small}
A category $\cC$ is finally (resp. cofinally) small
if and only if it admits a small final (resp. cofinal)
subcategory.
\end{corollary}
\begin{proof} Indeed, the condition is obviously sufficient.
Conversely, suppose that $\phi:I\to\cC$ is a final
(resp. cofinal) functor, with $I$ small; the datum
$\cG:=(\phi(\Ob(I)),\phi(\rMorph(I)))$ is a subgraph
of the category $\cC$ (see \eqref{subsec_graphs}), and
the inclusion morphism into $\cC$ extends to a functor
$\phi':I'\to\cC$, where $I'$ is the path category of
$\cG$ (proposition \ref{prop_graph-to-cat}). It is then
easily seen that the graph $(\phi(\Ob(I')),\phi'(\rMorph(I')))$
is actually a final (resp. cofinal) small subcategory
of $\cC$ (details left to the reader).
\end{proof}

The path category of a graph $\cG$ is a sort of
free category generated by $\cG$. We explain next
how to define equivalence relations on such a path
category, and how to form the corresponding quotient
categories. To this aim, let $\cG$ be any graph; we
shall call a {\em commutativity condition} on $\cG$
any pair of paths of $\cG$ with the same sources and
targets. A {\em conditional graph} is a pair $(\cG,C)$
consisting of a graph $\cG$ and a set $C$ of commutativity
conditions on $\cG$. A {\em morphism} $F:(\cG,C)\to(\cG',C')$
of conditional graphs is a morphism $F:\cG\to\cG'$ of
graphs, which induces a map $C\to C'$ :
$(p_1,p_2)\mapsto(P(F)p_1,P(F)p_2)$. With this composition
rule, the conditional graphs form also a category, and
we say that a conditional graph is {\em small} if its
underlying graph is small. We denote
$$
\sU\tdu\mathbf{CondGraph}
$$
the category of all small conditional graphs, and as
usual we drop the subscript $\sU$, unless ambiguities
are likely to arise from this omission.

\sset\subsubsection{}
We may attach to any small category $\cC$ a natural
set of commutativity conditions; namely, one takes
the set $C_\cC$ consisting of :
\begin{itemize}
\item
all pairs of paths
$((A_0,f_1,\dots,A_n),(B_0,g_1,\dots,B_m))$ in $\cC$
such that $A_0=B_0$, $A_n=B_m$, and
$f_n\circ\cdots\circ f_1=g_m\circ\cdots\circ g_0$
\item
the pairs $((A,\one_A,A),(A))$, for $A$ ranging over
all the objects of $\cC$.
\end{itemize}
Clearly, every functor $F:\cC\to\cC'$ induces a
map $C_\cC\to C_{\cC'}$, whence a forgetful functor
\set\begin{equation}\label{eq_cat-to-condgraph}
\bCat\to\mathbf{CondGraph}
\qquad
\cC\mapsto(\cC,C_\cC).
\end{equation}

\begin{proposition}\label{prop_cat-to-condgraph}
The forgetful functor \eqref{eq_cat-to-condgraph}
admits a left adjoint.
\end{proposition}
\begin{proof} Given a conditional graph $(\cG,C)$,
denote by $\cP$ the path category of $\cG$, and by
$\cQ$ the category such that
\begin{itemize}
\item
$\Ob(\cQ)=V(\cG)$
\item
for every $A,B\in\Ob(\cQ)$, the set $\Hom_\cQ(A,B)$
consists of all pairs $(p_1,p_2)\in P(\cG)\times P(\cG)$
such that the source (resp. the target) of both $p_1$ and
$p_2$ is $A$ (resp. is $B$)
\item
the composition law of $\cQ$ is given by the rule :
$(p_1,p_2)\circ(p'_1,p'_2):=(p_1\circ p'_1,p_2\circ p'_2)$
for every pair of composable morphisms
$(p_1,p_2),(p'_1,p'_2)\in\rMorph(\cQ)$
\end{itemize}
and notice that $C\subset\rMorph(\cQ)$. Consider the family
$S_C$ of all subcategories $\cS$ of $\cQ$ such that
\begin{itemize}
\item
$\Ob(\cS)=\Ob(\cQ)$ and $C\subset\rMorph(\cS)$
\item
$\rMorph(\cS)$ is an equivalence relation on
$P(\cG)=\rMorph(\cP)$.
\end{itemize}
Clearly $\cQ\in S_C$, so $S_C\neq\emptyset$, and $S_C$
admits a smallest element, namely the subcategory
$\cR$ whose set of objects is $\Ob(\cQ)$ and whose set of
morphisms is $\bigcap_{\cS\in S_C}\rMorph(\cS)$. We may then
form a category
$$
\cP/C
$$
whose set of objects is $\Ob(\cP)=V(\cG)$, and such that
$$
\Hom_{\cP/C}(A,B):=\Hom_\cP(A,B)/\!\sim_C
$$
where $\sim_C$ is the equivalence relation induced by $\cR$
on $\Hom_\cP(A,B)$. The composition law for morphisms is the
unique one such that the resulting pair of maps
$$
\Ob(\cP)\to\Ob(\cP/C)
\qquad
\rMorph(\cP)\to\rMorph(\cP/C)
$$
(which is the identity map on objects, and the natural
projection on morphisms) is a functor
$$
\cP\to\cP/C.
$$
To check that this law is well defined, say that
$p,q:A\to A'$ and $p',q':A'\to A''$ are two pairs of
morphisms in $\cP$, such that $p\sim_Cq$ and $p'\sim_Cq'$;
since $\cR$ is a subcategory, it follows easily that
$p'\circ p\sim_C q'\circ p\sim_Cq'\circ q$, whence
the claim. Lastly, let $\cC$ be any small category, and
$F:(\cG,C)\to(\cC,C_\cC)$ a given morphism of conditional
graphs; by proposition \ref{prop_graph-to-cat}, the morphism
$F$ extends to a functor $G:\cP\to\cC$. Let $\cS\subset\cQ$
be the subcategory such that $\Ob(\cS)=\Ob(\cQ)$ and whose
morphisms are all the pairs of paths $(p,p')$ such that
$(P(F)p,P(F)p')\in C_\cC$; by construction, $\cS\in S_C$,
so $\cS$ contains $\cR$; it follows that $G$ factors
uniquely through $\cP/C$, whence the proposition.
\end{proof}

\begin{theorem}\label{th_localize-cats}
Let $\cC$ be any category and $\Sigma\subset\rMorph(\cC)$
any subset. Then there exist a category $\cC[\Sigma^{-1}]$
and a {\em localization functor}
$$
L:\cC\to\cC[\Sigma^{-1}]
$$
with the following properties :

{\em(i)}\ \
$Lf$ is an isomorphism in $\cC[\Sigma^{-1}]$, for every
$f\in\Sigma$.

{\em(ii)}\ \
For every category $\cB$ and any functor $G:\cC\to\cB$
such that $Gf$ is an isomorphism in $\cB$ for every
$f\in\Sigma$, there exists a unique functor
$G':\cC[\Sigma^{-1}]\to\cB$ such that $G=G'\circ L$.
\end{theorem}
\begin{proof} Pick any universe $\sU'$ containing $\sU$,
and such that $\cC$ is $\sU'$-small; after replacing $\sU$
by $\sU'$, we may assume that $\cC$ is small. We consider
the graph $\cG$ such that
\begin{itemize}
\item
$V(\cG)=\Ob(\cC)$
\item
$\cG(A,B)$ is the disjoint union of $\Hom_\cC(A,B)$
and $\Hom_\cC(B,A)\cap\Sigma$, for every $A,B\in V(\cG)$.
\end{itemize}
For every $f\in\Hom_\cC(B,A)\cap\Sigma$, we write
$f^{-1}:A\to B$ for the corresponding arrow of $\cG$.
Next, we endow $\cG$ with the set $\Theta$ consisting
of the following commutativity conditions :
\begin{itemize}
\item
For every $A\in\Ob(\cC)$, the pair $((A,\one_A,A),(A))$
\item
for every $A,B,C\in\Ob(\cC)$, every $f\in\Hom_\cC(A,B)$
and every $g\in\Hom_\cC(B,C)$, the pair
$((A,f,B,g,C),(A,g\circ f,C))$
\item
the pairs $((A,f^{-1},B,f,A),(A,\one_A,A))$ and
$((B,f,A,f^{-1},B),(B,\one_B,B))$, for every
$A,B\in\Ob(\cC)$ and every $f\in\Hom_\cC(B,A)\cap\Sigma$.
\end{itemize}
By proposition \ref{prop_cat-to-condgraph}, the left adjoint
of \eqref{eq_cat-to-condgraph} maps the conditional graph
$(\cG,\Theta)$ to a category which we denote $\cC[\Sigma^{-1}]$.
Now, suppose that $G:\cC\to\cB$ is a functor such that
$Gf$ is an isomorphism in $\cB$, for every $f\in\Sigma$.
We define a morphism of graphs $H:\cG\to\cB$ by the
rule :
\begin{itemize}
\item
$H(A):=GA$ for every $A\in V(\cG)$
\item
$H(f):=Gf$ for every $A,B\in V(\cG)$ and every
$f\in\Hom_\cC(A,B)$
\item
$H(f):=(Gf)^{-1}$ for every $A,B\in V(\cG)$ and every
$f\in\Hom_\cC(B,A)\cap\Sigma$.
\end{itemize}
It is easily seen that $H$ induces a morphism of
conditional graphs
$$
(\cG,\Theta)\to(\cB,C_\cB)
$$
which, under the adjunction of proposition
\ref{prop_cat-to-condgraph}, corresponds to a unique
functor $\cC[\Sigma^{-1}]\to\cB$ with the sought properties.
\end{proof}

\begin{remark}\label{rem_local-cat}
(i)\ \
Clearly conditions (i) and (ii) of theorem
\ref{th_localize-cats} characterize $\cC[\Sigma^{-1}]$ up
to isomorphism of categories. 
We call a {\em localization of\/ $\cC$ relative to $\Sigma$}
the datum of a functor $G:\cC\to\cD$ and an equivalence
of categories $H:\cC[\Sigma^{-1}]\to\cD$ such that $H\circ L=G$,
where $L:\cC\to\cC[\Sigma^{-1}]$ is the localization functor
of theorem \ref{th_localize-cats}.

(ii)\ \
For any subset $\Sigma\subset\rMorph(\cC)$, let us set
$\Sigma^o:=\{s^o~|~s\in\Sigma\}\subset\rMorph(\cC^o)$
(notation of \eqref{subsec_opposing}). In view of (i),
we see that the localization functor $\cC^o\to\cC^o[\Sigma^{o-1}]$
extends uniquely to an isomorphism of categories :
$$
\cC[\Sigma^{-1}]^o\isom\cC^o[\Sigma^{o-1}].
$$
\end{remark}

\begin{corollary}\label{cor_of-localization}
In the situation of theorem {\em\ref{th_localize-cats}}, the
localization functor $L$ induces a fully faithful functor
$$
\bFun(L,\cB):\bFun(\cA[\Sigma^{-1}],\cB)\to\bFun(\cA,\cB)
\qquad
\text{for every category $\cB$}.
$$
\end{corollary}
\begin{proof} Let $F,G:\cA[\Sigma^{-1}]\to\cB$ be two functors;
recall that a natural transformation $\alpha:F\Rightarrow G$
is the same as the datum of a functor
$\tilde\alpha:\cA[\Sigma^{-1}]\to\sMorph(\cB)$ with
$\ss\circ\tilde\alpha=F$ and $\st\circ\tilde\alpha=G$ (see
\eqref{subsec_Morph-cat}). Notice also that a morphism
$(g,g')$ of $\sMorph(\cB)$ is invertible if and only if both
$\ss(g,g'):=g$ and $\st(g,g'):=g'$ are invertible in $\cB$.
Combining with the universal property of $F$ from theorem
\ref{th_localize-cats}, we conclude that the datum of $\alpha$
is equivalent to that of a functor $\tilde\beta:\cA\to\sMorph(\cB)$
such that $\ss\circ\tilde\beta=F\circ L$ and
$\st\circ\tilde\beta=G\circ L$; but the
latter is the same as a natural transformation
$F\circ L\Rightarrow G\circ L$, whence the contention.
\end{proof}

\begin{example}\label{ex_invert-all-with-final}
Let $\cC$ be a category that admits either a final or an
initial object, and let $\Sigma$ be the set of all morphisms
of $\cC$. Then $\cC[\Sigma^{-1}]$ is (isomorphic to) the
category $\overline\cC$ whose set of objects is $\Ob(\cC)$
and such that for every $X,Y\in\Ob(\cC)$ the set
$\Hom_{\overline\cC}(X,Y)$ contains exactly one element.
Indeed, there exists a unique functor $p:\cC\to\bar\cC$
which is the identity on objects; since every morphism of
$\bar\cC$ is an isomorphism, $p$ factors uniquely through
a functor $\cC[\Sigma^{-1}]\to\bar\cC$ that is still the
identity on objects. On the other hand, say that $X_0$
is a final object for $\cC$, and for every $X\in\Ob(\cC)$,
let $t_X:X\to X_0$ be the unique morphism; denote by
$[t_X]$ the class of $t_X$ in $\cC[\Sigma^{-1}]$. Then
we have also a functor $q:\bar\cC\to\cC[\Sigma^{-1}]$
that is identity on objects, and sends every morphism
$X\to Y$ of $\bar\cC$ to the morphism
$\tau_{XY}:=[t_Y]^{-1}\circ[t_X]:X\to Y$ of $\cC[\Sigma^{-1}]$.
Clearly $p\circ q=\one_{\bar\cC}$, and to conclude, it suffices
to check that if $f:Y\to X$ is any morphisms of $\cC$, then
the class $[f]:X\to Y$ of $f$ in $\cC[\Sigma^{-1}]$ coincides
with $\tau_{XY}$. But the latter assertion is clear, since
$t_Y\circ f=t_X$ in $\cC$. One argues likewise in case $\cC$
admits an initial object.
\end{example}

\begin{proposition}\label{prop_local-and-adj-pairs}
Let $\cC\!,\!\cD$ be two categories,
$F\!:\!\cC\to\cD,G\!:\!\cD\to\cC$ two functors, and set
$$
\Sigma_F:=\{f\in\rMorph(\cC)~|~\text{$Ff$ is an isomorphism}\}.
$$
Suppose that $G$ is fully faithful, and is either right or
left adjoint to $F$. Then $F$ factors uniquely through the
localization $L:\cC\to\cC[\Sigma_F^{-1}]$ and an equivalence
$$
F':\cC[\Sigma_F^{-1}]\isom\cD.
$$
\end{proposition}
\begin{proof} To ease notation, set $\cC':=\cC[\Sigma_F^{-1}]$,
and notice that $\Sigma_{F^o}=\Sigma_F^o$ (notation of remark
\ref{rem_local-cat}(ii)), whence a natural isomorphism
of categories :
$$
\cC'^o\isom\cC^o[\Sigma^{-1}_{F^o}]
$$
that identifies $L^o:\cC^o\to\cC'^o$ with the localization
$\cC^o\to\cC^o[\Sigma^{-1}_{F^o}]$; moreover, recall that $(F,G)$
is an adjoint pair if and only if the same holds for the
pair $(G^o,F^o)$ (remark \ref{rem_opposite-Fun}(iv)). Thus,
we may assume that $G$ is right adjoint to $F$. Now, the
existence and uniqueness of $F'$ is clear from the universal
property of the localization; it remains to check that $F'$
is an equivalence. To this aim, set as well
$G':=L\circ G:\cD\to\cC'$, and let $(\eta,\eps)$ be the unit
and counit of an adjunction for the pair $(F,G)$. Since
$F=F'\circ L$, we get $F'\circ G'=F\circ G$, hence $\eps$
is also a natural transformation $F'G'\Rightarrow\one_\cD$.
On the other hand, by corollary \ref{cor_of-localization}
there exists a unique natural transformation
$$
\eta':\one_{\cC'}\Rightarrow G'F'
\qquad\text{such that}\qquad
\eta'*L=L*\eta.
$$
Using the triangular identities of \eqref{subsec_adj-pair},
we compute :
$$
(G'*\eps)\odot(\eta'*G')=(LG*\eps)\odot(L*\eta*G)=
L*((G*\eps)\odot(\eta*G))=L*\one_G=\one_{G'}.
$$
Likewise, let us show that $(\eps*F')\odot(F'*\eta')=\one_{F'}$; to
this aim, again by virtue of corollary \ref{cor_of-localization},
it suffices to check that
$$
((\eps*F')\odot(F'*\eta'))*L=\one_{F'}*L=\one_F.
$$
However, the left-hand side equals
$(\eps*F'L)\odot(F'*\eta'*L)=(\eps*L)\odot(F'L*\eta)$, so
the sought identity follows again from \eqref{subsec_adj-pair}.
By proposition \ref{prop_triangular-identities}(i), it
follows that $(F',G')$ is an adjoint pair of functors,
with unit $\eta'$ and counit $\eps$. Lastly, since $G$
is fully faithful, $\eps$ is an isomorphism of functors
(proposition \ref{prop_fullfaith-adjts}(i)); invoking
again the triangular identities \eqref{subsec_adj-pair},
we deduce that $F*\eta$ is an isomorphism of functors,
hence $\eta_X\in\Sigma_F$ for every $X\in\Ob(\cC)$. Since
$\Ob(\cC)=\Ob(\cC')$, it follows easily that $\eta'$ is
an isomorphism of functors, whence the proposition.
\end{proof}

\begin{definition}\label{def_right-calculus}
Let $\cC$ be a category, $\Sigma\subset\rMorph(\cC)$
a set of morphisms.
\begin{enumerate}
\item
We say that $\Sigma$ {\em admits a right calculus of fractions}
if the following conditions hold :
\begin{itemize}
\item[(CF1)]
$\one_A\in\Sigma$, for every $A\in\Ob(\cC)$.
\item[(CF2)]
For every $s:A\to B$ and $t:B\to C$ with $s,t\in\Sigma$,
we have $t\circ s\in\Sigma$ as well.
\item[(CF3)]
For every morphism $f:A\to B$ in $\cC$ and every $s:C\to B$
in $\Sigma$, there exist $g:D\to C$ in $\cC$ and
$t:D\to A$ in $\Sigma$ such that $f\circ t=s\circ g$.
\item[(CF4)]
If $f,g:A\to B$ are any two morphisms in $\cC$ such that
$s\circ f=s\circ g$ for some $s:B\to C$ in $\Sigma$, then
there exists $t:D\to A$ in $\Sigma$ such that
$f\circ t=g\circ t$.
\end{itemize}
\item
We say that $\Sigma$ {\em admits a left calculus of fractions}
if the subset $\Sigma^o:=\{s^o~|~s\in\Sigma\}$ admits a right
calculus of fractions (notation of remark \ref{rem_local-cat}(ii)).
\item
For every $A\in\Ob(\cC)$, let $\Sigma_A$ be the full
subcategory of $\cC/A$ whose objects are the elements
of $\Sigma$ with target equal to $A$ (notation of
\eqref{subsec_slice-cat}). We say that $\Sigma$ is
{\em right cofinally small}, if $\Sigma_A$ is cofinally
small for every $A\in\Ob(\cC)$ (see definition
\ref{def_MacL-cofinal}(iv)).
\end{enumerate}
\end{definition}

\sset\subsubsection{}\label{subsec_right-calculus}
Let now $\cC$ and $\Sigma$ be as in definition
\ref{def_right-calculus}, and suppose that $\cC$ has
small $\Hom$-sets, and $\Sigma$ admits a right calculus
of fractions.  For every $A,B\in\Ob(\cC)$, let us consider
the functor
$$
H_{A,B}:\Sigma_A^o\to\Set
\qquad
(s:I\to A)\mapsto\{(s,f)~|~f\in\Hom_\cC(I,B)\}
$$
where, for every morphism $h/A:s\to t$ in $\Sigma_A$,
the mapping $H_{A,B}(h/A)$ is given by the rule
$$
H_{A,B}(t)\to H_{A,B}(s)
\quad :\quad
(t,f)\mapsto(s,f\circ h).
$$
Furthermore, we define a mapping
$$
c^{A,B}_s:H_{A,B}(s)\to\Hom_{\cC[\Sigma^{-1}]}(A,B)
\qquad
(s,f)\mapsto f\circ s^{-1}
\qquad
\text{for every $s\in\Ob(\Sigma_A)$}.
$$

\begin{proposition}\label{prop_calculus-frac}
In the situation of \eqref{subsec_right-calculus},
the following holds :
\begin{enumerate}
\item
For every $A\in\Ob(\cC)$, the category $\Sigma^o_A$
is filtered.
\item
For every $A,B\in\Ob(\cC)$, the rule :
$$
s\mapsto c^{A,B}_s
\qquad
\text{for every $s\in\Ob(\Sigma_A)$}
$$
defines a universal cocone with basis $H_{A,B}$ and vertex
$\Hom_{\cC[\Sigma^{-1}]}(A,B)$.
\item
If\/ $\Sigma$ is right cofinally small, The category
$\cC[\Sigma^{-1}]$ has small $\Hom$-sets.
\end{enumerate}
\end{proposition}
\begin{proof}(i): First, $\Sigma^o_A$ is not empty, due
to (CF1). Next, say that $s,s'\in\Ob(\Sigma_A)$, and
denote by $I$ and $I'$ the sources of $s$ and respectively
$s'$; due to (CF3), we may find $I''\in\Ob(\cC)$, a morphism
$f:I''\to I$ in $\cC$, and an element $g:I''\to I'$ of
$\Sigma$ such that $s\circ f=s'\circ g$. Then, (CF2) says
that $s'\circ g$ lies in $\Sigma$, and therefore it defines
an object $s''$ of $\Sigma_A$, with morphisms $f:s\to s''$
and $g:s'\to s''$ in $\Sigma_A^o$. Lastly, say that
$f,f':s\to s'$ are any two morphisms in $\Sigma^o_A$, and
denote by $I'$ the source of $s'$. By (CF4), we may find
$g:I''\to I'$ in $\Sigma$ such that $f\circ g=f'\circ g$;
clearly $s'\circ g$ defines a morphism $h:s'\to s''$ in
$\Sigma^o_A$ such that $h\circ f=h\circ f'$. Now the assertion
follows from remark \ref{rem_cofinal}(i).

(ii): Let $s:I\to A$ and $s':I'\to A$ be two objects
of $\Sigma_A$, and $g:I'\to I$ a morphism $s\to s'$ in
$\Sigma^o_A$. Given $f:I\to B$, set $f':=f\circ g$; we
notice :

\begin{claim}\label{cl_inverting-s}
$f\circ s^{-1}=f'\circ s'{}^{-1}$ in $\cC[\Sigma^{-1}]$.
\end{claim}
\begin{pfclaim} Due to (CF3), we may find $I''\in\Ob(\cC)$,
an element $t:I''\to I$ of $\Sigma$ and a morphism
$t':I''\to I'$ such that $s\circ t=s'\circ t'$.
We compute :
$$
s\circ g\circ t'=s'\circ t'=s\circ t
$$
whence, by (CF4), an element $h:I'''\to I''$
of $\Sigma$, such that $g\circ t'\circ h=t\circ h$.
It follows that $g\circ t'=t$ in $\cC[\Sigma^{-1}]$
and therefore
$$
f\circ t=f\circ g\circ t'=f'\circ t'
\qquad
\text{in $\cC[\Sigma^{-1}]$}.
$$
Consequently :
$$
f\circ s^{-1}=f'\circ t'\circ t^{-1}\circ s^{-1}=
f'\circ s'{}^{-1}\circ s\circ t\circ t^{-1}\circ s^{-1}=
f'\circ s'{}^{-1}
\qquad
\text{in $\cC[\Sigma^{-1}]$}
$$
as stated.
\end{pfclaim}

Claim \ref{cl_inverting-s} says already that $c_\bullet^{A,B}$
is a well defined cocone. To check the universality property,
we define a new category $\cD$ with $\Ob(\cD)=\Ob(\cC)$ and
such that $\Hom_\cD(A,B)$ is the colimit of the functor
$H_{A,B}$, for every $A,B\in\Ob(\cC)$. From (i) and example
\ref{ex_complete-cats}(i) we see that this is the set of
equivalence classes $[s,f]$ of pairs $(s,f)\in H_{A,B}(s)$
with $s$ ranging over all objects of $\Sigma_A$; two such
pairs $(s,f)$ and $(s',f')$ are equivalent, if there exist
$t\in\Ob(\Sigma_A)$ and morphisms $h:t\to s$, $h':t\to s'$
in $\Sigma_A$ such that $f\circ h=f'\circ h'$.
The composition of morphisms in $\cD$ is defined
as follows. Let $A,B,C\in\Ob(\cC)$ be any three objects,
$s:I\to A$, $t:J\to B$ any elements of $\Sigma$, and
$(s,f)\in H_{A,B}(s)$, $(t,g)\in H_{B,C}(t)$; by (CF3)
we find $D\in\Ob(\cC)$, a morphism $f':D\to J$ in $\cC$
and an element $t':D\to I$ of $\Sigma$ such that
$t\circ f'=f\circ t'$. We define
\set\begin{equation}\label{eq_define-composition}
(t,g)\circ(s,f):=[s\circ t',g\circ f']\in\Hom_\cD(A,C).
\end{equation}
Let us check that this class does not depend on the choice
of $D$, $f'$ and $t'$. Indeed, suppose that $D'\in\Ob(\cC)$
is another object and $t'':D'\to I$ in $\Sigma$, $f'':D''\to J$
in $\cC$ satisfy the condition $t\circ f''=f\circ t''$;
by (CF3) we may then find $D''\in\Ob(\cC)$, $h:D''\to D$
in $\Sigma$ and $h':D''\to D'$ in $\cC$ such that
$t'\circ h=t''\circ h'$. We compute :
$$
t\circ f'\circ h=f\circ t'\circ h=f\circ t''\circ h'=
t\circ f''\circ h'.
$$
By (CF4), it follows that there exist $D'''\in\Ob(\cC)$
and an element $u:D'''\to D''$ of $\Sigma$ such that
$f'\circ h\circ u=f''\circ h'\circ u$. By (CF2), we may then
replace $D''$ by $D'''$ and $h$, $h'$ with $h\circ u$ and
respectively $h'\circ u$, after which we may assume as well
that $f''\circ h'=f'\circ h$. Finally, we find :
$$
[s\circ t',g\circ f']=[s\circ t'\circ h,g\circ f'\circ h]
=[s\circ t''\circ h',g\circ f''\circ h']=[s\circ t'',g\circ f'']
$$
as required. It also follows that $[s\circ t',g\circ f']$
depends only on $[s,f]$ and $[t,g]$; indeed, say that
$u:I'\to I$ is any element of $\Sigma$, and set
$(s',f''):=(s\circ u,f\circ u)$. Suppose also that
$f''\circ t''=t\circ f'$ for some $t'':D\to I'$ in
$\Sigma$; then
$(s'\circ t'',g\circ f')=(s\circ u\circ t'',g\circ f')$,
and the foregoing shows that the latter pair is
equivalent to $(s\circ t',g\circ f')$, which shows
the independence on the chosen representative for
$[s,f]$; similarly one checks the independence on the
representative for $[t,g]$.

We next check the associativity of the composition law
thus obtained. To this aim, say that $A,B,C,D\in\Ob(\cC)$
are any four objects, $s:I\to A$, $t:J\to B$, $u:K\to C$
any three elements of $\Sigma$, and $f:I\to B$, $g:J\to C$,
$h:K\to D$ any three morphisms of $\cC$. Choose
$E\in\Ob(\cC)$ with an element $t':E\to I$ of $\Sigma$
and a morphism $f':E\to J$ of $\cC$ such that
$f\circ t'=t\circ f'$, and therefore
$$
[t,g]\circ[s,f]=[s\circ t',g\circ f'].
$$
Then, choose $F\in\Ob(\cC)$ with $u':F\to J$ in $\Sigma$
and $g':F\to K$ in $\cC$ such that $g\circ u'=u\circ g'$,
and therefore
$$
[u,h]\circ[t,g]=[t\circ u',h\circ g'].
$$
By (CF3) we may find $G\in\Ob(\cC)$, $u'':G\to E$ in $\Sigma$
and $f'':G\to F$ in $\cC$ such that $u'\circ f''=f'\circ u''$.
By (CF2), $t'\circ u''$ lies in $\Sigma$, and therefore it is
easily seen that the pair $(s\circ t'\circ u'',h\circ g'\circ f'')$
represents both $[t\circ u',h\circ g']\circ[s,f]$ and
$[u,h]\circ[s\circ t',g\circ f']$, whence the assertion.
Let us also remark that -- due to (CF1) -- for every
$A\in\Ob(\cC)$, the class $[\one_A,\one_A]$ is the
identity endomorphism of $A$, seen as an object of $\cD$.
This completes the construction of the category $\cD$.
Next, we claim that the system of cocones
$c^{\bullet\bullet}_\bullet$ yields a functor
$$
F:\cD\to\cC[\Sigma^{-1}]
$$
whose map on objects is the identity mapping of $\Ob(\cC)$,
and such that
$$
F[s,f]:=f\circ s^{-1}
\qquad
\text{for every morphism $[s,f]$ in $\cD$}.
$$
Indeed, claim \ref{cl_inverting-s} implies that the
foregoing rule yields a well defined map
$$
\Hom_\cD(A,B)\to\Hom_{\cC[\Sigma^{-1}]}(A,B)
\qquad
\text{for every $A,B\in\Ob(\cC)$}.
$$
It is also clear that $F\one_A=\one_A$ for every
$A\in\Ob(\cC)$, hence it remains only to check that
$F$ respects the composition laws of the two categories;
however, say that $s,f,t,g$ are as in
\eqref{eq_define-composition}; then we have
$$
(g\circ f')\circ(s\circ t')^{-1}=
(g\circ t^{-1}\circ f\circ t')\circ t'{}^{-1}\circ s^{-1}=
(g\circ t^{-1})\circ(f\circ s^{-1})
\qquad
\text{in $\cC[\Sigma^{-1}]$}
$$
whence the contention. Likewise, we have a functor
$$
G:\cC\to\cD
$$
which is the identity mapping on objects, and such that
$Gf:=[\one_A,f]$ for every $A,B\in\Ob(\cC)$ and every
morphism $f:A\to B$ in $\cC$. Morever, the identities
$$
[\one_I,s]\circ[s,\one_I]=[\one_A,\one_A]
\qquad
[s,\one_I]\circ[\one_I,s]=[\one_I,\one_I]
\qquad
\text{for every $s:I\to A$ in $\Sigma$}
$$
show that $Fs$ is invertible in $\cD$, for every such $s$.
By theorem \ref{th_localize-cats}, it follows that
$G$ factors through a unique functor
$G':\cC[\Sigma^{-1}]\to\cD$, and a simple inspection
shows that $F\circ G'$ is the identity endofunctor
of $\cC[\Sigma^{-1}]$. Likewise, it is easily seen
that $G'\circ F=\one_\cD$, so $F$ establishes an
isomorphism of categories, and (ii) follows.

(iii) is an immediate consequence of (ii).
\end{proof}

\begin{remark}\label{rem_left-calc-fract}
Proposition \ref{prop_calculus-frac}(ii) means that
$c^{A,B}_\bullet$ induces a natural identification
$$
\colim_{\Sigma^o_A}H_{A,B}\isom\Hom_{\cC[\Sigma^{-1}]}(A,B)
\qquad
\text{for every $A,B\in\Ob(\cC)$}.
$$
The dual of this assertion provides a corresponding
computation of the morphisms in
$\cC[\Sigma^{-1}]$ in terms of left fractions. Namely,
suppose that $\Sigma$ admits a left calculus of fractions;
for every $A\in\Ob(\cC)$ we let $\Sigma^A$ be the full
subcategory of $A/\cC$ whose objects are the elements
of $\Sigma$ with source equal to $A$, and we say that
$\Sigma$ is {\em left cofinally small} if $\Sigma^A$
is cofinally small for every such $A$. For every other
object $B$ of $\cC$ we consider the functor
$$
H_{B,A}:\Sigma^A\to\Set
\qquad
(s:A\to I)\mapsto\{(f,s)~|~f\in\Hom_\cC(B,I)\}
$$
where, for every morphism $h:s\to t$ in $\Sigma^A$,
the mapping $H_{B,A}(h)$ is given by the rule
$$
H_{B,A}(s)\to H_{B,A}(t)
\quad :\quad
(f,s)\mapsto(h\circ f,t).
$$
Then $\Sigma^A$ is a filtered category, and if $\Sigma$
is left cofinally small, $\cC[\Sigma^{-1}]$ has small
$\Hom$-sets; more precisely, we have a natural
identification :
$$
\colim_{\Sigma^A}H_{B,A}\isom\Hom_{\cC[\Sigma^{-1}]}(B,A)
\qquad
(f,s)\mapsto s^{-1}\circ f.
$$
\end{remark}

\sset\subsubsection{}\label{subsec_full-subcat-fract}
Let $\cA$ be a category, $\cA_0$ a full subcategory of
$\cA$, $\Sigma\subset\Ob(\cA)$ any subset, and set
$$
\Sigma_0:=\Sigma\cap\rMorph(\cA_0).
$$
Clearly, the inclusion functor $i:\cA_0\to\cA$ extends
uniquely to a functor
$$
i[\Sigma_0^{-1}]:\cA_0[\Sigma_0^{-1}]\to\cA[\Sigma^{-1}].
$$

\begin{proposition}\label{prop_full-subcat-fract}
In the situation of \eqref{subsec_full-subcat-fract},
suppose that :
\begin{enumerate}
\alphaenu
\item
For every $A\in\Ob(\cA)$, $A_0\in\Ob(\cA_0)$ and every
morphism $f:A\to iA_0$ that lies in $\Sigma$, there exist
$A'_0\in\Ob(\cA_0)$ and a morphism $g:iA'_0\to A$ such that
$f\circ g\in\Sigma$.
\item
$\Sigma$ admits a right calculus of fractions.
\end{enumerate}
Then the set $\Sigma_0$ admits a right calculus of
fractions in $\cA_0$, and the functor
$i[\Sigma_0^{-1}]$ is fully faithful.
\end{proposition}
\begin{proof} As usual, after replacing our universe $\sU$
by a larger one, we may assume that $\cA$ is small. Let us
check that $\Sigma_0$ admits a right calculus of fractions.
(CF1) and (CF2) are obviously fulfilled by $\Sigma_0$. Next,
suppose that we have $A,B,C\in\Ob(\cA_0)$, $D\in\Ob(\cA)$
and a commutative diagram in $\cA$
\set\begin{equation}\label{eq_CF3}
{\diagram
D \ar[r]^-g \ar[d]_t & iC \ar[d]^{is} \\
iA \ar[r]^-{if} & iB
\enddiagram}
\qquad
\text{with $s,t\in\Sigma$}.
\end{equation}
By (a), it follows that there exists a morphism
$u:iD_0\to D$ for some $D_0\in\Ob(\cA_0)$, such that
$t':=t\circ u:iD_0\to iA$ lies in $\Sigma_0$,
$g':=g\circ u:iD_0\to iC$ lies in $\Hom_{\cA_0}(D_0,C)$,
and clearly $f\circ t'=s\circ g'$, whence (CF3).
Lastly, suppose that $A,B\in\Ob(\cA_0)$,
$D\in\Ob(\cA)$, and we have morphisms $f,g:A\to B$ and
$t:D\to iA$ such that $i(f)\circ t=i(g)\circ t$ and with
$t\in\Sigma$. Pick $u:iD_0\to D$ as in the foregoing,
and let again $t':=t\circ u$; then  $f\circ t'=g\circ t'$,
and we have just seen that $t'\in\Sigma_0$, whence (CF4).

Next, let us check that $i[\Sigma^{-1}_0]$ is fully faithful.
In light of propositions \ref{prop_calculus-frac}(ii)
and \ref{prop_MacL-cofinal}, it suffices to show that the
inclusion functor
\set\begin{equation}\label{eq_many-sigmas}
\Sigma^o_{0,A}\to\Sigma^o_{iA}
\end{equation}
is cofinal for every $A\in\Ob(\cA_0)$ (notation of
\eqref{subsec_right-calculus}). To this aim, in view
of proposition \ref{prop_calculus-frac}(i), it
suffices to check that the functor \eqref{eq_many-sigmas}
fulfills conditions (a) and (b) of lemma
\ref{lem_filtered-final}(i). However, condition (a)
of {\em loc.cit.} translates directly as the assumption
(a) of the proposition. Next, consider any pair of morphisms
$f,g:iD_0\to D$ and any element $s:D\to iA$ in $\Sigma$ such
that $s\circ f=s\circ g$. By (CF4), there exists $t:E\to iD_0$
in $\Sigma$ such that $f\circ t=g\circ t$, and then
assumption (a) yields a morphism $u:iE_0\to E$ such
that $t':=t\circ u$ lies in $\Sigma_0$, and
$f\circ i(t')=g\circ i(t')$, which shows that condition
(b) of lemma \ref{lem_filtered-final}(iv) holds as well.
\end{proof}

\begin{remark}\label{rem_full-subcat-fract}
Taking into account remarks \ref{rem_local-cat}(ii)
and \ref{rem_left-calc-fract}, we see that the dual of
proposition \ref{prop_full-subcat-fract} holds as well;
{\em i.e.} in the situation of \eqref{subsec_full-subcat-fract},
suppose that
\begin{enumerate}
\alphaenu
\item
For every $A\in\Ob(\cA)$, $A_0\in\Ob(\cA_0)$ and every
morphism $f:iA_0\to A$ that lies in $\Sigma$, there exist
$A'_0\in\Ob(\cA)$ and a morphism $g:A\to iA'_0$ such that
$g\circ f\in\Sigma$.
\item
$\Sigma$ admits a left calculus of fractions.
\end{enumerate}
Then the set $\Sigma_0$ admits a left calculus of
fractions in $\cA_0$, and the functor $i[\Sigma_0^{-1}]$
is fully faithful.
\end{remark}

\begin{proposition}\label{prop_rcf-and-cofilt-commas}
In the situation of \eqref{subsec_right-calculus},
let $F:\cC\to\cC[\Sigma^{-1}]$ be the localization functor.
The following holds :
\begin{enumerate}
\item
For every $X\in\Ob(\cC[\Sigma^{-1}])$ the category $X/F\cC$
is cofiltered.
\item
If\/ $\cC$ is small, the functor
$F_!:\cC^\wedge\to\cC[\Sigma^{-1}]^\wedge$ is exact.
\end{enumerate}
\end{proposition}
\begin{proof}(i): Let $((Y_i,f_i:X\to FY_i)~|~i=1,2)$ be a
pair of objects of $X/F\cC$. For $i=1,2$, we may write
$f_i=g_i\circ s_i^{-1}$ with some morphism $s_i:Z_i\to X$ in
$\Sigma$ and some morphism $g_i:Z_i\to Y_i$ in $\cC$. By
condition (CF3) of definition \ref{def_right-calculus} we
may then find a morphism $t_1:Z\to Z_1$ in $\Sigma$ and a
morphism $t_2:Z\to Z_2$ in $\cC$ such that
$s:=s_1\circ t_1=s_2\circ t_2$, and notice that $s\in\Sigma$,
by condition (CF2). For $i=1,2$, the composition $t_i\circ g_i$
then yields a morphism $(Z,s^{-1}:X\to FZ)\to(Y_i,f_i)$ in
$X/F\cC$. Next, let $h_1,h_2:(Y_1,f_1)\to(Y_2,f_2)$ be two
morphisms in $X/F\cC$; this means that $h_i:Y_1\to Y_2$ is a
morphism in $\cC$ for $i=1,2$, and $Fh_1\circ f_1=Fh_2\circ f_1=f_2$
in $\cC[\Sigma^{-1}]$. Thus, $F(h_1\circ g_1)=F(h_2\circ g_1)$,
and therefore there exists $(u:Z'\to Z_1)\in\Sigma$ such that
$h_1\circ g_1\circ u=h_2\circ g_1\circ u$ in $\cC$. Then the
composition $g:=g_1\circ u$ yields a morphism
$g:(Z',(s_1\circ u)^{-1}:X\to FZ')\to(Y_1,f_1)$ in $X/F\cC$
such that $h_1\circ g=h_2\circ g$. The assertion then follows
from remark \ref{rem_cofinal}(i).

(ii): The assertion follows from (i), arguing as in the proof
of corollary \ref{cor_lable}(i).
\end{proof}

\section{\texorpdfstring{$2$}{2}-Category theory}
In dealing with categories, the notion of equivalence is
much more central than the notion of isomorphism. On the other
hand, equivalence of categories is usually not preserved by the
standard categorical operations discussed thus far. For instance,
consider the following :

\begin{example} Let $\cC$ be the category with $\Ob(\cC)=\{a,b\}$,
and whose only morphisms are $\one_a$, $\one_b$ and $u:a\to b$,
$v:b\to a$. Then necessarily $u\circ v=\one_b$ and $v\circ u=\one_a$.
Let $\cC_a$ (resp. $\cC_b$) be the unique subcategory of $\cC$
with $\Ob(\cC_a)=\{a\}$ (resp. $\Ob(\cC_b)=\{b\}$). Clearly both
inclusion functors $\cC_a\to\cC\leftarrow\cC_b$ are equivalences.
However, $\cC_a\times_\cC\cC_b$ is the empty category; especially,
this fibre product is not equivalent to $\cC=\cC\times_\cC\cC$.
\end{example}

It is therefore natural to seek a new framework for the manipulation
of categories and functors ``up to equivalences'', and thus more
consonant with the very spirit of category theory. Precisely such
a framework is provided by the theory of $2$-categories, the subject
of this chapter.

\subsection{\texorpdfstring{$2$}{2}-Categories and pseudo-functors}
\label{sec_2Cats}
The category $\bCat$, together with the category structure on the
sets $\Fun(-,-)$ (as in \eqref{subsec_Godem-prod}), provides the
first example of a {\em $2$-category}. The latter is the datum of :
\begin{itemize}
\item
A set $\Ob(\cA)$, whose elements are called the
{\em objects of $\cA$}.
\item
For every $A,B\in\Ob(\cA)$, a category $\cA(A,B)$.
The objects of $\cA(A,B)$ are called {\em $1$-cells\/} or
{\em arrows}, and are designated by the usual arrow notation
$f:A\to B$.
Given $f,g\in\Ob(\cA(A,B))$, we shall write $f\Rightarrow g$
to denote a morphism from $f$ to $g$ in $\cA(A,B)$. Such morphisms
are called {\em $2$-cells}. The composition of $2$-cells
$\alpha:f\Rightarrow g$ and $\beta:g\Rightarrow h$ shall
be denoted by $\beta\odot\alpha:f\Rightarrow h$.
\item
For every $A,B,C\in\Ob(\cA)$, a {\em composition bifunctor} :
$$
c_{ABC}:\cA(A,B)\times\cA(B,C)\to\cA(A,C).
$$
Given $1$-cells $A\xrightarrow{f}B\xrightarrow{g}C$, we write
$g\circ f:=c_{ABC}(f,g)$. \\
Given two $2$-cells $\alpha:f\Rightarrow g$ and $\beta:h\Rightarrow k$,
respectively in $\cA(A,B)$ and $\cA(B,C)$, we use the notation
$$
\beta*\alpha:=c_{ABC}(\alpha,\beta):h\circ f\Rightarrow k\circ g.
$$
Also, if $h$ is any $1$-cell of $\cA(B,C)$, we usually write
$h*\alpha$ instead of $\one_h*\alpha$. Likewise, we set
$\beta*f:=\beta*\one_f$, for every $1$-cell $f$ in $\cA(A,B)$.
\item
For every element $A\in\Ob(\cA)$, a {\em unit functor} :
$$
u_A:\one\to\cA(A,A)
$$
where $\one:=(*,\one_*)$ is the final object of $\bCat$.
Hence $u_A$ is the datum of an object:
$$
\one_A\in\Ob(\cA(A,A))
$$
and its identity endomorphism, which we shall denote by
$i_A:\one_A\to\one_A$.
\end{itemize}
The bifunctors $c_{ABC}$ are required to satisfy an
{\em associativity axiom}, which says that the diagram:
$$
\xymatrix{
\cA(A,B)\times\cA(B,C)\times\cA(C,D)
\ar[rr]^-{\one\times c_{BCD}} \ar[d]_{c_{ABC}\times\one} & &
\cA(A,B)\times\cA(B,D) \ar[d]^{c_{ABD}} \\
\cA(A,C)\times\cA(C,D) \ar[rr]^-{c_{ACD}} & & \cA(A,D)
}$$
commutes for every $A,B,C,D\in\Ob(\cA)$. Likewise, the functor
$u_A$ is required to satisfy a {\em unit axiom}; namely,
the diagram :
$$
\xymatrix{
\one\times\cA(A,B) \ar[d]_{u_A\times\one_{\cA(A,B)}} &
\cA(A,B) \ar[l]_-\sim \ar[r]^-\sim \ddouble &
\cA(A,B)\times\one \ar[d]^{\one_{\cA(A,B)}\times u_B} \\
\cA(A,A)\times\cA(A,B) \ar[r]^-{c_{AAB}} &
\cA(A,B) & \cA(A,B)\times\cA(B,B) \ar[l]_-{c_{ABB}}
}$$
commutes for every $A,B\in\Ob(\cA)$.

\begin{remark}\label{rem_equiv-2-cat}
(i)\ \
In any $2$-category, consider a diagram with two $2$-cells :
$$
\xymatrix{ A \rtwocell^f_g{\alpha} &
B \rtwocell^{f'}_{g'}{\alpha'} & C.
}$$
In this situation, we have the following commutation
identity :
$$
(\alpha'*g)\odot(f'*\alpha)=\alpha'*\alpha=
(g'*\alpha)\odot(\alpha'*f)
$$
which is sometimes useful, to perform certain verifications.
To see this, notice that $\alpha'*\alpha=
(\alpha'\odot\one_{f'})*(\one_g\odot\alpha)$, whence the
first stated identity, by functoriality of the composition
for $2$-cells. Likewise, we get
$\alpha'*\alpha=(\one_{g'}\odot\alpha')*(\alpha\odot\one_f)$,
whence the second identity.

(ii)\ \
As already announced, the category $\bCat$ becomes naturally
a $2$-category, by letting
$$
\bCat(\cA,\cB):=\bFun(\cA,\cB)
\qquad
\text{for every $\cA,\cB\in\Ob(\bCat)$}
$$
with composition bifunctor defined as in
\eqref{subsec_Godem-prod}. Indeed, by spelling out the
definitions, we see that the functoriality of the
Godement product boils down precisely to the identity
\eqref{eq_Godement-functor}.

On the other hand, from every $2$-category $\cA$ we obtain
a category, simply by forgetting the $2$-cells : the set of
objects of this {\em underlying category of $\cA$} is
$\Ob(\cA)$, and the morphisms are the $1$-cells of $\cA$,
with composition law induced by the composition bifunctor
of $\cA$.

(iii)\ \
We shall say that a $2$-category $\cA$
{\em has $\sU$-small $\Hom$-categories} if $\cA(A,B)$
is a $\sU$-small category for every $A,B\in\Ob(\cA)$.
We shall say that $\cA$ is {\em $\sU$-small} if it has
$\sU$-small $\Hom$-categories and $\Ob(\cA)$ is $\sU$-small.
As usual, the universe $\sU$ will be dropped from the
notation, unless there is a danger of ambiguities. For
instance, the $2$-category $\bCat$ has small $\Hom$-categories.
\end{remark}

\sset\subsubsection{}\label{subsec_opp-2-cats}
If $\cA$ is a $2$-category, we have three distinct
constructions of {\em opposite $2$-categories\/}
associated with $\cA$. Namely we have :
\begin{itemize}
\item
The $2$-category $\cA^o$ such that
$\Ob(\cA^o):=\{A^o~|~A\in\Ob(\cA)\}$ (cp.
\eqref{subsec_opposing}) and with
$$
\cA^o(A^o,B^o):=\cA(B,A)
\qquad
\text{for every $A^o,B^o\in\Ob(\cA^o)$}.
$$
For any $1$-cell $f:B\to A$ in $\cA$, we denote
$f^o:A^o\to B^o$ the corresponding $1$-cell in $\cA^o$.
The composition bifunctor $c^o_{\bullet\bullet\bullet}$
and unit functors $u^o_\bullet$ of $\cA^o$ are given
by the rules :
$$
c^o_{A^oB^oC^o}:=c_{CBA}
\qquad
\one_{A^o}:=\one^o_A
\qquad
\text{for every $A^o,B^o,C^o\in\Ob(\cA^o)$}.
$$
Notice that to any $2$-cell $\alpha:f\Rightarrow g$
in $\cA$ there corresponds a $2$-cell
$\alpha^o:f^o\Rightarrow g^o$ in $\cA^o$. The associativity
and unit axioms for $\cA^o$ follow formally from
the corresponding properties for $\cA$.
\item
The $2$-category ${}^o\!\cA$ such that
$\Ob({}^o\!\cA):=\{{}^o\!\!A~|~A\in\Ob(\cA)\}$ and with
$$
{}^o\!\cA({}^o\!A,{}^o\!B):=\cA(A,B)^o
\qquad
\text{for every ${}^o\!\!A,{}^o\!B\in\Ob({}^o\!\cA)$}.
$$
For any $1$-cell $f:A\to B$ in $\cA$, we denote
${}^o\!f:{}^o\!A\to {}^o\!B$ the corresponding $1$-cell
in ${}^o\!\cA$. The composition bifunctor
${}^o\!c_{\bullet\bullet\bullet}$
and unit functors ${}^o\!u_\bullet$ of ${}^o\!\cA$ are given
by the rules :
$$
{}^o\!c_{{}^o\!A{}^oB{}^oC}:=c^o_{ABC}
\qquad
\one_{{}^o\!A}:={}^o\one_A
\qquad
\text{for every ${}^o\!\!A,{}^o\!B,{}^o\!C\in\Ob({}^o\!\cA)$}.
$$
Then to any $2$-cell $\alpha:f\Rightarrow g$ in $\cA$ there
corresponds a $2$-cell
${}^o\alpha:{}^og\Rightarrow{}^o\!f$ in ${}^o\cA$. The
associativity and unit axioms for ${}^o\!\cA$ are also
immediately deduced from those of $\cA$.
\item
Lastly, we may combine the two foregoing constructions
to get the $2$-category ${}^o\!\cA^o$, which inverts the
direction of both $1$-cells and $2$-cells.
\end{itemize}

\begin{definition}\label{def_adjoint-1-cells}
Let $\cA$ be any $2$-category, and $f:A\to B$ any $1$-cell
of $\cA$.

(i)\ \
A $2$-cell $\beta:f\Rightarrow f'$ of $\cA$ is
{\em invertible}, if it is an isomorphism in
the category $\cA(A,B)$.

(ii)\ \
We say that the $1$-cell $g:B\to A$ of $\cA$ is
{\em right adjoint} to $f$, if there exist $2$-cells
$\eta:\one_A\Rightarrow g\circ f$ and
$\eps:f\circ g\Rightarrow\one_B$ fulfilling the
{\em triangular identities}
$$
(g*\eps)\odot(\eta*g)=\one_g
\qquad
(\eps*f)\odot(f*\eta)=\one_f.
$$
(see \cite[Def.7.1.2]{Bor}). In this case, we also say
that $g$ is {\em right adjoint} to $f$, and that $(f,g)$
is an {\em adjoint pair of\/ $1$-cells}; moreover, we say
that $(\eta,\eps)$ is an {\em adjunction} for $(f,g)$.
Then $\eta$ is called the {\em unit} and $\eps$ the
{\em counit} of the adjunction $(\eta,\eps)$.

(iii)\ \
We say that a $1$-cell $f:A\to B$ of $\cA$ is an
{\em equivalence\/} if there exist $1$-cells $g,h:B\to A$
and invertible $2$-cells $\one_A\Rightarrow g\circ f$
and $f\circ h\Rightarrow\one_B$. In this case, we also
say that $g$ and $h$ are {\em quasi-inverse} $1$-cells
for $f$.
\end{definition}

\begin{remark}\label{rem_compose-equiv}
(i)\ \
Let $f,g$ and $h$ be as in definition
\ref{def_adjoint-1-cells}(iii), so we have invertible
$2$-cells $\eta:\one_A\Rightarrow g\circ f$ and
$\eps:f\circ h\Rightarrow\one_B$. Then we have the
invertible $2$-cell
$\beta:=(g*\eps)\odot(h*\eta):h\Rightarrow g$, whence
the invertible $2$-cell
$\eps\circ(f*\beta^{-1}):f\circ g\Rightarrow\one_B$.
In other words, we may as well assume that $g=h$ in
definition \ref{def_adjoint-1-cells}(iii).

(ii)\ \
A composition of equivalences is an equivalence. Indeed,
let $f_1:A_1\to A_2$ and $f_2:A_2\to A_3$ be two equivalences,
$g_1:A_2\to A_1$ (resp. $g_2:A_3\to A_2$) a quasi-inverse for
$f_1$ (resp. for $f_2$) so that we have invertible $2$-cells
$\alpha_i:\one_{A_i}\Rightarrow g_i\circ f_i$ and
$\beta_i:f_i\circ g_i\Rightarrow\one_{A_{i+1}}$ for $i=1,2$.
Then $g:=g_1\circ g_2$ is a quasi-inverse for $f:=f_2\circ f_1$,
since we have the invertible $2$-cells
$(g_1*\alpha_2*f_1)\odot\alpha_1:\one_{A_1}\Rightarrow g\circ f$
and $\beta_2\odot(f_2*\beta_1*g_2):f\circ g\Rightarrow\one_{A_3}$.

(iii)\ \
Let $f,f':A\to B$ be two $1$-cells of $\cA$, and $\beta:f\isom f'$
an invertible $2$-cell. Then $f$ is an equivalence if and only if the
same holds for $f'$. Indeed, suppose we have a $1$-cell $g:B\to A$
with an invertible $2$-cells $\alpha:\one_A\isom g\circ f$ and
$\alpha':f\circ g\isom\one_B$; then we get the invertible $2$-cells
$(g*\beta)\odot\alpha:\one_A\isom g\circ f'$ and
$\alpha'\odot(\beta^{-1}*g):f'\circ g\isom\one_B$, so $f'$ is an
equivalence.

(iv)\ \
Likewise, if $(f,g)$ is an adjoint pair of $1$-cells of $\cA$
and we have invertible $2$-cells $\alpha:g\isom g'$ and
$\beta:f\isom f'$, then $(f',g')$ is an adjoint pair of
$1$-cells. Indeed, by assumption we have $2$-cells
$\eta:\one_A\Rightarrow g\circ f$ and $\eps:f\circ g\Rightarrow\one_B$
fulfilling the triangular identities of definition
\ref{def_adjoint-1-cells}(ii); it suffices to check that
$\eta':=(\alpha*\beta)\odot\eta$ and
$\eps':=\eps\odot(\beta*\alpha)^{-1}$ are a unit and a counit
for the pair $(f',g')$. However, we have :
$$
\begin{aligned}
(g'*\eps')\odot(\eta'*g')
&\,=(g'*(\eps\odot(\beta*\alpha)^{-1}))
\odot(((\alpha*\beta)\odot\eta)*g') \\
&\,=(g'*\eps)\odot(g'*\beta*\alpha)^{-1}\odot(\alpha*\beta*g')
\odot(\eta*g') \\
&\,=(g'*\eps)\odot(\alpha*f*\alpha^{-1})\odot(\eta*g') \\
&\,=(g'*\eps)\odot(\alpha*f*g)\odot(g*f*\alpha^{-1})\odot(\eta*g') \\
&\,=\alpha\odot(g*\eps)\odot(\eta*g)\odot\alpha^{-1} \\
&\,=\alpha\odot\alpha^{-1}=\one_{g'}
\end{aligned}
$$
where the third, fourth and fifith equalities follow from remark
\ref{rem_equiv-2-cat}(i). Likewise one checks the other required
triangular identity.

(v)\ \
Furthermore, if $f:A\to B$ is a $1$-cell in $\cA$ and
$g,g':B\to A$ are both right adjoint to $f$, then we have
an invertible $2$-cell $\alpha:g\isom g'$. Indeed, let
$(\eta,\eps)$ (resp. $(\eta',\eps')$) be an adjunction
for the pair $(f,g)$ (resp. for the pair $(f,g')$). We set
$$
\alpha:=(g'*\eps)\odot(\eta'*g):g\Rightarrow g'
\qquad
\beta:=(g*\eps')\odot(\eta*g'):g'\Rightarrow g.
$$
Applying repeatedly remark \ref{rem_equiv-2-cat}(i), we compute :
$$
\begin{aligned}
\beta\odot\alpha
&\,=(g*\eps')\odot(g*f*g'*\eps)\odot(\eta*g'*f*g)\odot(\eta'*g) \\
&\,=(g*\eps)\odot(g*\eps'*f*g)\odot(g*f*\eta'*g)\odot(\eta*g) \\
&\,=(g*\eps)\odot(g*((\eps'*f)\odot(f*\eta'))*g)\odot(\eta*g) \\
&\,=(g*\eps)\odot(\eta*g)=\one_g.
\end{aligned}
$$
Likewise one checks that $\alpha\odot\beta=\one_{g'}$, whence
the contention.
\end{remark}

\sset\subsubsection{}\label{subsec_square-algebra}
Let $\cA$ be any $2$-category; an {\em oriented square
in $\cA$} is a (not necessarily commutative) diagram of
the type :
\set\begin{equation}\label{eq_2-cell}
{\diagram A \ar[r]^-f \ar[d]_g
\drtwocell\omit{_\ \alpha} & B \ar[d]^h \\
A' \ar[r]_-{f'} & B'
\enddiagram}\end{equation}
where $f,f'g,h$ are any four $1$-cells and
$\alpha:h\circ f\Rightarrow f'\circ g$ is any $2$-cell
of $\cA$. If $\alpha$ is an invertible $2$-cell, we also
say that the square \eqref{eq_2-cell} is
{\em essentially commutative}. We say that $\alpha$
{\em orients} the square \eqref{eq_2-cell}, and often
we refer to such a square just by naming its orienting
$2$-cell; when dealing with complex diagrams, this
shorthand may lead to unacceptable ambiguities, but
in these cases usually one may resolve such ambiguities
just by specifying two of the opposite sides of the square,
for which we shall employ a fractional notation : thus,
the square \eqref{eq_2-cell} shall be denoted, depending
on the context, in either of the following three manners :
$$
\alpha
\qquad
\frac{\alpha}{f|f'}
\qquad
\frac{\alpha}{g|h}.
$$
One basic operation consists in combining two squares that
share one edge, to obtain a new square, essentially by
omitting the common edge of the squares. Namely, suppose
we have two squares $\alpha$ and $\beta$ as in the diagram :
\set\begin{equation}\label{eq_star-comp}
{\diagram A \ar[r]^-g \ar[d]_f
\drtwocell\omit{^\alpha\ } & A' \ar[d]|{f'} \ar[r]^-k
\drtwocell\omit{^\beta\ } & A'' \ar[d]^{f''} \\
B \ar[r]_-h & B' \ar[r]_-i & B''
\enddiagram}
\end{equation}
Then we say that $\alpha$ and $\beta$ are {\em vertically
composable}, and we define the square
$$
\beta\boxvert\alpha
\qquad :\qquad
{\diagram A \ar[r]^-{k\circ g} \ar[d]_f
\drtwocell\omit{^\gamma\ } & A'' \ar[d]^{f''} \\
B \ar[r]_{i\circ h} & B''
\enddiagram}
\qquad
\text{where\ \ $\gamma:=(\beta*g)\odot(i*\alpha)$}.
$$
In this operation, we have joined to the square $\alpha$
of \eqref{eq_2-cell} another square $\beta$ that shares
the edge $f'$. Clearly, if $\beta'$ is another square
that shares with $\alpha$ the edge $f$, the same operation
can be performed to obtain the square $\alpha\boxvert\beta'$
which omits the edge $f$. On the other hand, the rule to
join to $\alpha$ a square that shares one of the two
remaining edges $g$ or $h$, is slightly different. Namely,
suppose we have the diagram of two squares :
$$
\xymatrix{ A \ar[r]^-f \ar[d]_g
\drtwocell\omit{_\ \alpha} & B \ar[d]|h 
\ar[r]^-i
\drtwocell\omit{_\ \gamma} & C \ar[d]^k \\
A' \ar[r]_-{f'} & B' \ar[r]_{i'} & C'.
}$$
Then we say that $\alpha$ and $\gamma$ are {\em horizontally
composable}, and we set
$$
\frac{\gamma}{i|i'}\boxminus\frac{\alpha}{f|f'}:=
\frac{(i'*\alpha)\odot(\gamma*f)}{i\circ f|i'\circ f'}.
$$
However, notice that the second join operation is turned
into the first one, when we replace $\cA$ by its opposite
$2$-categories; namely, we have the identities :
\set\begin{equation}\label{eq_change-orientation}
(\gamma\boxminus\alpha)^o=\alpha^o\boxvert\gamma^o
\qquad
{}^o(\gamma\boxminus\alpha)={}^o\gamma\boxvert{}^o\alpha
\end{equation}
The following proposition establishes the two basic
rules of the algebra of oriented squares.

\begin{proposition}\label{prop_square-algebra}
With the notation of \eqref{subsec_square-algebra},
the following holds :
\begin{enumerate}
\item
The join operation is associative, {\em i.e.} for
any diagram of three squares
$$
\xymatrix{ A \ar[r]^-g \ar[d]_f
\drtwocell\omit{^\alpha\ } & A' \ar[d]|{f'} \ar[r]^-k
\drtwocell\omit{^\beta\ } & A'' \ar[d]|{f''} \ar[r]^-l
\drtwocell\omit{^\gamma\ } & A''' \ar[d]^{f'''} \\
B \ar[r]_-h & B' \ar[r]_-i & B'' \ar[r]_-m & B'''
}$$
we have
$$
\gamma\boxvert(\beta\boxvert\alpha)=
(\gamma\boxvert\beta)\boxvert\alpha.
$$
\item
For every (not necessarily commutative) diagram
$$
\xymatrix{
A \ar[r]^-f \ar[d]_{h}
\drtwocell\omit{_\alpha\ } & B \ar[d]|k \ar[r]^-g
\drtwocell\omit{_\beta\ } & C \ar[d]^i \\
A' \ar[r]|{f'} \ar[d]_{h'}
\drtwocell\omit{_\alpha'\ } & B' \ar[r]|{g'} \ar[d]|{k'}
\drtwocell\omit{_\beta'\ } & C' \ar[d]^{i'} \\
A'' \ar[r]_-{f''} & B'' \ar[r]_-{g''} & C''.
}$$
we have
$$
(\beta'\boxvert\beta)\boxminus(\alpha'\boxvert\alpha)=
(\beta'\boxminus\alpha')\boxvert(\beta\boxminus\alpha).
$$
\end{enumerate}
\end{proposition}
\begin{proof}(i): We compute :
$$
\begin{aligned}
\gamma\boxvert(\beta\boxvert\alpha)=\, &
\gamma\boxvert((\beta*g)\odot(i*\alpha)) \\
=\, & ((\gamma*(k\circ g))\odot(m*((\beta*g)\odot(i*\alpha))) \\
=\, &
((\gamma*(k\circ g))\odot(m*(\beta*g))\odot(m*(i*\alpha)) \\
=\, &
(((\gamma*k)\odot(m*\beta))*g)\odot((m\circ i)*\alpha) \\
=\, & (\gamma\boxvert\beta)\boxvert\alpha.
\end{aligned}
$$

(ii): We compute :
$$
\begin{aligned}
(\beta'\boxvert\beta)\boxminus(\alpha'\boxvert\alpha)
=\, & (g''*((\alpha'*h)\odot(k'*\alpha)))\odot
(((\beta'*k)\odot(i'*\beta))*f) \\
=\, & (g''*\alpha'*h)\odot(g''*k'*\alpha)\odot(\beta'*k*f)
\odot(i'*\beta*f) \\
=\, & (g''*\alpha'*h)\odot(\beta'*f*h)\odot(i'*g'*\alpha)
\odot(i'*\beta*f) \\
=\, & ((\beta'\boxminus\alpha')*h)\odot
(i'*(\beta\boxminus\alpha)) \\
=\, & (\beta'\boxminus\alpha')\boxvert(\beta\boxminus\alpha)
\end{aligned}
$$
where the fourth equality follows from the commutation
rule of remark \ref{rem_equiv-2-cat}(i).
\end{proof}

\begin{remark}\label{rem_square-algebra}
Proposition \ref{subsec_square-algebra}(i) establishes
explicitly only the associativity of the join operation
$\boxvert$, but in light of \eqref{eq_change-orientation}
we also get immediately the associativity of $\boxminus$.
\end{remark}

\begin{example}\label{ex_2-arrows}
(i)\ \
Let $\cA$ be any $2$-category. We may construct a new
$2$-category
$$
2\tdu\sMorph(\cA)
$$
which is the $2$-categorical counterpart of the category
of arrows (see \eqref{subsec_Morph-cat}). Namely : 
\begin{itemize}
\item
the objects of $2\tdu\sMorph(\cA)$ are the $1$-cells of
$\cA$.
\item
If $f:A\to B$ and $f':A'\to B'$ are any two $1$-cells of
$\cA$, the arrows $f\to f'$ in $2\tdu\sMorph(\cA)$ are
all the oriented squares \eqref{eq_2-cell}.
\item
For $1$-cells
$(g_1,h_1,\alpha_1),(g_2,h_2,\alpha_2):f\to f'$, the
$2$-cells $(g_1,h_1,\alpha_1)\Rightarrow(g_2,h_2,\alpha_2)$
are all the pairs
$$
(\beta:g_1\Rightarrow g_2,\gamma:h_1\Rightarrow h_2)
\qquad\text{such that}\qquad
\frac{\alpha_1}{h_1|g_1}\boxminus\frac{\beta}{g_1|g_2}=
\frac{\gamma}{h_1|h_2}\boxminus\frac{\alpha_2}{h_2|g_2}.
$$
\item
For $f$ as in the foregoing, the unit functor
$\one\to 2\tdu\sMorph(f,f)$ is given by the $2$-cell
\set\begin{equation}\label{eq_unit-in-morph}
{\diagram A \ar[r]^-f \ar[d]_{\one_A}
\drtwocell\omit{_\ \ \one_f} & B \ar[d]^{\one_B} \\
A \ar[r]_-f & B
\enddiagram}
\end{equation}
and its identity endomorphism $(i_A,i_B)$.
\end{itemize}
For $f$ and $f'$ as in the foregoing, the composition law
in $2\tdu\sMorph(\cA)(f,f')$ is given by the rule
$$
(\beta',\gamma')\odot(\beta,\gamma):=
(\beta'\odot\beta,\gamma'\odot\gamma)
$$
for every composable pair $((\beta,\gamma),(\beta',\gamma'))$
of $2$-cells in $2\tdu\sMorph(\cA)$. The composition
bifunctor is given by the composition law for squares;
{\em i.e.}, given $1$-cells $f\to f$' and $f'\to f''$
in $2\tdu\sMorph(\cA)$ as in \eqref{eq_star-comp}, we set
$$
(\beta,k,i)\circ(\alpha,g,h):=
(\beta\boxvert\alpha,k\circ g,i\circ h).
$$
Lastly, given $2$-cells $(\beta,\gamma)$ in
$2\tdu\sMorph(\cA)(f,f')$ and $(\beta',\gamma')$ in
$2\tdu\sMorph(\cA)(f',f'')$ we set
$$
(\beta',\gamma')*(\beta,\gamma):=(\beta'*\beta,\gamma'*\gamma).
$$
The functoriality and associativity properties of the
composition laws for $2$-cells are then immediate from
the definitions. The associativity of the composition
law for $1$-cells holds by proposition
\ref{prop_square-algebra}(i).

(ii)\ \
Just as in \eqref{subsec_slice-cat}, for every $X\in\Ob(\cA)$
we may consider the $2$-subcategories $\cA/X$ and $X/\cA$ of
$2\tdu\sMorph(\cA)$. Explicitly, the objects of $\cA/X$ are
the $1$-cells $A\to X$ of $\cA$. Given two such objects $f,f'$,
the $1$-cells $f\to f'$ of $\cA/X$ are all the diagrams in
$\cA$ of the type
$$
\xymatrix{ A \ar[r]^-f \ar[d]_g
\drtwocell\omit{_\ \alpha} & X \ddouble \\
A' \ar[r]_-{f'} & X
}$$
and for $1$-cells $(g_1,\alpha_1),(g_2,\alpha_2):f\to f'$,
the $2$-cells of $\cA/X$ are the $2$-cells
$\beta:g_1\Rightarrow g_2$ of $\cA$ such that
$(f'*\beta)\odot\alpha_1=\alpha_2$. A similar description applies
to $X/\cA$.

(iii)\ \
To every $1$-cell $h:X\to Y$ in $\cA$ we attach functors
$$
h^Z_*:\cA(Z,X)\to\cA(Z,Y)
\qquad
\text{for every $Z\in\Ob(\cA)$}
$$
by ruling that $h^Z_*(k):=h\circ k$ and $h^Z_*(\beta):=h*\beta$
for every $1$-cell $k:Z\to X$ of $\cA$ and every morphism
$\beta:k\Rightarrow k'$ in $\cA(Z,X)$. Let $\nu:h\Rightarrow h'$
be any $2$-cell in $\cA$; we claim that the rule
$$
k\mapsto(\nu^Z_{*k}:=\nu*k:h^Z_*(k)\Rightarrow h'^Z_*(k))
\qquad
\text{for every $(k:Z\to X)\in\Ob(\cA(Z,X))$}
$$
defines a natural transformation
$$
\nu^Z_*:h^Z_*\Rightarrow h'^Z_*
\qquad
\text{for every $Z\in\Ob(\cA)$}.
$$
Indeed, let $k,k':Z\to X$ be any two objects of $\cA(Z,X)$,
and $\mu:k\Rightarrow k'$ a $2$-cell; we need to show that
the resulting diagram commutes in $\cA$ :
$$
\xymatrix{ h^Z_*(k) \ar@{=>}[r]^-{\nu*k} \ar@{=>}[d]_{h*\mu} &
h'^Z_*(k) \ar@{=>}[d]^{h'*\mu} \\
h^Z_*(k') \ar@{=>}[r]^-{\nu*k'} & h'^Z_*(k')
}$$
But this follows directly from remark \ref{rem_equiv-2-cat}(i).
Likewise, we may attach to $f$ the functors
$$
h_Z^*:\cA(Y,Z)\to\cA(X,Z)
\qquad
k\mapsto k\circ h
\qquad
\beta\mapsto\beta*h
\qquad
\text{for every $Z\in\Ob(\cA)$}
$$
and then $\nu$ induces as well a natural transformation
$$
\nu^*_Z:h^*_Z\Rightarrow h'^*_Z
\qquad
k\mapsto k*\nu
\qquad
\text{for every $Z\in\Ob(\cA)$}.
$$
We apply these constructions to prove the following result, which
generalizes propositions \ref{prop_triangular-identities}(ii)
and \ref{prop_fullfaith-adjts}(i) to arbitrary $2$-categories :
\end{example}

\begin{lemma}\label{lem_long-time}
Let $\cA$ be a $2$-category, $f:A\to B$, $g:B\to A$
two $1$-cells in $\cA$. We have :
\begin{enumerate}
\item
Suppose that $f$ is an equivalence, and $g$ a quasi-inverse
for $f$. Then, for every invertible $2$-cell
$\eta:\one_A\Rightarrow g\circ f$ there exists a unique
invertible $2$-cell $\eps:f\circ g\Rightarrow\one_B$ such
that $(\eta,\eps)$ is an adjunction for the pair $(f,g)$.
\item
Suppose that $(f,g)$ is an adjoint pair, and that $(\eta,\eps)$
and $(\eta',\eps')$ are two adjunctions for $(f,g)$. Then
$\eta=\eta'$ if and only if $\eps=\eps'$.
\item
The following conditions are equivalent :
\begin{enumerate}
\alphaenu
\item
$f$ is an equivalence.
\item
For every $Z\in\Ob(\cA)$ the functor $f^*_Z:\cA(B,Z)\to\cA(A,Z)$
is an equivalence.
\item
For every $Z\in\Ob(\cA)$ the functor $f^Z_*:\cA(Z,A)\to\cA(Z,B)$
is an equivalence.
\end{enumerate}
\end{enumerate}
\end{lemma}
\begin{proof}(i): By assumption, there exists an invertible
$2$-cell $\eps':f\circ g\Rightarrow\one_B$. Now, since
$\eta$ and $\eps'$ are invertible $2$-cells, it is clear
that $\eta^Z_*$ and $\eps'^Z_*$ are isomorphisms of
functors for every $Z\in\Ob(\cA)$, so $f^Z_*$ and $g^Z_*$
are equivalences of categories (proposition
\ref{prop_fullfaith-adjts}(i)); therefore there exists
a unique adjunction for the pair $(f^Z_*,g^Z_*)$ whose unit
is $\eta^Z_*$ (claim \ref{cl_get-adjunction}) and we denote
by $\eps^Z_*$ the unique counit for this adjunction.
Then also $\eps^Z_*$ is an isomorphism of functors
(proposition \ref{prop_fullfaith-adjts}(iii)), and we
define the invertible $2$-cell
$\eps:=\eps^B_{\one_B}:f\circ g\Rightarrow\one_B$. Lastly,
the triangular identities \eqref{subsec_adj-pair} for the
pair $(\eta^Z_*,\eps^Z_*)$ immediately imply the
corresponding identities for $(\eta,\eps)$.

(ii): As in the foregoing, we obtain for every $Z\in\Ob(\cA)$
an adjoint pair of functors $(f^Z_*,g^Z_*)$, and two adjunctions
$(\eta^Z_*,\eps^Z_*)$, $(\eta'^Z_*,\eps'^Z_*)$ for $(f^Z_*,g^Z_*)$.
Then proposition \ref{prop_triangular-identities}(ii) implies
that $\eta^Z_*=\eta'^Z_*$ if and only if $\eps^Z_*=\eps'^Z_*$.
The assertion is an immediate consequence.

(iii.a)$\Rightarrow$(iii.b),(iii.c): Indeed, let $g:B\to A$ be
a quasi-inverse for $f$; it is easily seen that the functor
$g^*_Z$ is a quasi-inverse for $f^*_Z$ and $g^Z_*$ is a
quasi-inverse for $f^Z_*$.

(iii.b)$\Rightarrow$(iii.a): Since $f^*_A$ is an equivalence,
we find a $1$-cell $g:B\to A$ of $\cA$ with an isomorphism
$g\circ f= f^*_A(g)\isom\one_A$. There follows an isomorphism
$f^*_B(f\circ g)=f\circ g\circ f\isom f=f^*_B(\one_B)$,
and since $f^*_B$ is an equivalence, we deduce an isomorphism
$f\circ g\isom\one_B$, whence the contention. Likewise one shows
that (iii.c)$\Rightarrow$(iii.a), or alternatively one can deduce
this implication from the foregoing one, by considering the
opposite $2$-categories.
\end{proof}

\begin{definition}\label{def_pseudo-fun}
(\cite[Def.7.5.1]{Bor}) Let $\cA$ and $\cB$ be two $2$-categories.

(i)\ \
A {\em pseudo-functor\/} $F:\cA\to\cB$ is the datum of :
\begin{itemize}
\item
For every $A\in\Ob(\cA)$, an object $FA\in\Ob(\cB)$.
\item
For every $A,B\in\Ob(\cA)$, a functor :
$$
F_{AB}:\cA(A,B)\to\cB(FA,FB).
$$
We shall often omit the subscript, and write only $Ff$ instead of
$F_{AB}f:FA\to FB$, for a $1$-cell $f:A\to B$, and likewise for
$2$-cells.
\item
For every $A,B,C\in\Ob(\cA)$, a natural isomorphism $\gamma_{ABC}$
between two functors $\cA(A,B)\times\cA(B,C)\to\cB(FA,FC)$ as indicated
by the (not necessarily commutative) diagram :
\set\begin{equation}\label{eq_coh-constraint-one}
{\diagram
\cA(A,B)\times\cA(B,C) \ar[rrr]^-{c_{ABC}} \ar[d]_{F_{AB}\times F_{BC}}
& \drtwocell\omit{^\gamma_{ABC}\ \ \ \ \ }& &
\cA(A,C) \ar[d]^{F_{AC}} \\
\cB(FA,FB)\times\cB(FB,FC) \ar[rrr]_-{c_{FA,FB,FC}} & & & \cB(FA,FC).
\enddiagram}
\end{equation}
To ease notation, for every $(f,g)\in\cA(A,B)\times\cA(B,C)$, we shall
write $\gamma_{f,g}$ instead of $(\gamma_{ABC})_{(f,g)}:
c_{FA,FB,FC}(F_{AB}f,F_{BC}g)\Rightarrow F_{AC}(c_{ABC}(f,g))$.
\item
For every $A\in\Ob(\cA)$, a natural isomorphism $\delta_{\!A}$
between functors $\one\to\cB(FA,FA)$, as indicated by
the (not necessarily commutative) diagram :
\set\begin{equation}\label{eq_coh-constraint-two}
{\diagram
\one \ar[r]^-{u_A} \ddouble \drtwocell\omit{^\delta_{\!A}\ \ } &
\cA(A,A) \ar[d]^{F_{AA}} \\
\one \ar[r]_-{u_{FA}} & \cB(FA,FA).
\enddiagram}
\end{equation}
\end{itemize}
The system $(\delta_\bullet,\gamma_{\bullet\bullet\bullet})$
is called the {\em coherence constraint\/} for $F$.
This datum is required to satisfy:
\begin{itemize}
\item
A {\em composition axiom}, which says that the diagram
$$
\xymatrix{
Fh\circ Fg\circ Ff \ar@{=>}[rr]^{Fh*\gamma_{f,g}}
\ar@{=>}[d]_{\gamma_{g,h}*Ff} & &
Fh\circ F(g\circ f) \ar@{=>}[d]^{\gamma_{g\circ f,h}} \\
F(h\circ g)\circ Ff \ar@{=>}[rr]^{\gamma_{f,h\circ g}} & &
F(h\circ g\circ f)
}$$
commutes for every sequence of arrows
$A\xrightarrow{f}B\xrightarrow{g}C\xrightarrow{h}D$ in $\cA$.
We will write
$$
\gamma_{f,g,h}:Fh\circ Fg\circ Ff\Rightarrow F(h\circ g\circ f)
$$
for the common composition of these two pairs of cells.
\item
A {\em unit axiom}, which says that the diagrams :
$$
\xymatrix{ Ff\circ\one_{FA} \ar@{=>}[d]_{\one_{Ff}}
\ar@{=>}[rr]^-{Ff*\delta_{\!A}} & &
Ff\circ F\one_A \ar@{=>}[d]^{\gamma_{\one_{\!A},f}} &
\one_{FB}\circ Ff \ar@{=>}[d]_{\one_{Ff}}
\ar@{=>}[rr]^-{\delta_{\!B}*Ff} & &
F(\one_B)\circ Ff \ar@{=>}[d]^{\gamma_{f,\one_{\!B}}} \\
Ff \ar@{=>}[rr]^-{\one_{Ff}} & & F(f\circ\one_{\!A}) &
Ff \ar@{=>}[rr]^-{\one_{Ff}} & & F(\one_B\circ f)
}$$
commute for every arrow $f:A\to B$ (where, to ease
notation, we have written $\delta_A$ instead of
$(\delta_A)_*:\one_{FA}\Rightarrow F_{AA}\one_A$,
and likewise for $\delta_B$).
\end{itemize}

(ii)\ \
A pseudo-functor $F$ as in (i) is called {\em strict} if
\eqref{eq_coh-constraint-one} and \eqref{eq_coh-constraint-two}
both commute, and $\gamma_{ABC}$ and $\delta_A$ are the
identity natural transformations for every $A,B,C\in\Ob(\cA)$.
\end{definition}

\begin{remark}\label{rem_pseudo-funct}
With the notation of definition \ref{def_pseudo-fun}, we have :

(i)\ \
The functoriality of $F_{AB}$, for every $A,B\in\Ob(\cA)$,
comes down to the identities
$$
F(\one_f)=\one_{Ff}
\qquad
F(\beta\odot\alpha)=F\beta\odot F\alpha
\qquad
\text{for every $1$-cell $f$ and $2$-cells $\alpha,\beta$
of $\cA(A,B)$}.
$$

(ii)\ \
Likewise, the naturality of $\gamma_{ABC}$ boils down
to the commutativity of the diagram
$$
\xymatrix{ Fg\circ Ff \ar@{=>}[rr]^-{\gamma_{f,g}}
\ar@{=>}[d]_{F\beta*F\alpha} & &
F(g\circ f) \ar@{=>}[d]^{F(\beta*\alpha)} \\
Fg'\circ Ff' \ar@{=>}[rr]^-{\gamma_{f',g'}} & & F(g'\circ f')
}$$
for every $1$-cells $f,f':A\to B$, $g,g':B\to C$ and
$2$-cells $\alpha:f\Rightarrow f'$, $\beta:g\Rightarrow g'$.

(iii)\ \
If $h:C\to D$ and $h':C'\to D'$ are any two other $1$-cells,
and $\mu:h\Rightarrow h'$ is any $2$-cell, it follows easily
that the diagram
$$
\xymatrix{ Fh\circ Fg\circ Ff
\ar@{=>}[rr]^-{\gamma_{f,g,h}} \ar@{=>}[d]_{F\mu*F\beta*F\alpha} & &
F(h\circ g\circ f) \ar@{=>}[d]^{F(\mu*\beta*\alpha)} \\
Fh'\circ Fg'\circ Ff' \ar@{=>}[rr]^-{\gamma_{f',g',h'}} & &
F(h'\circ g'\circ f')
}$$
commutes : details left to the reader.

(iv)\ \
The pseudo-functor $F$ induces {\em opposite pseudo-functors}
$$
F^o:\cA^o\to\cB^o
\qquad\text{and}\qquad
{}^o\!F:{}^o\!\cA\to{}^o\!\cB
$$
(as well as the functor ${}^o\!F^o$, by combining the
two operations). Namely, $F^o$ is given by the rules
$$
F^o(A^o):=(FA)^o
\qquad
F^o(f^o):=(Ff)^o
\qquad
F^o(\beta^o):=(F\beta)^o
$$
for every $A\in\Ob(\cA)$, every $1$-cell $f$ and every
$2$-cell $\beta$ of $\cA$. The coherence constraint of
$F^o$ is the pair $(\gamma^o,\delta^o)$ such that
$$
\gamma^o_{g^o,f^o}:=(\gamma_{f,g})^o
\qquad\text{and}\qquad
\delta^o_{A^o}:=(\delta_A)^o
$$
for every composable pair of $1$-cells $f,g$ of $\cA$,
and every $A\in\Ob(\cA)$. The composition and unit axioms
for $F^o$ follow immediately from the same for $F$. Likewise,
${}^o\!F$ is given by the rules
$$
{}^o\!F({}^o\!A):={}^o\!(FA)
\qquad
{}^o\!F({}^o\!f):={}^o\!(Ff)
\qquad
{}^o\!F({}^o\!\beta):={}^o\!(F\beta)
$$
for every $A,f,\beta$ as in the foregoing. The coherence
constraint for ${}^o\!F$ is the pair $({}^o\!\gamma,{}^o\!\delta)$
with
$$
{}^o\!\gamma_{{}^o\!f,{}^o\!g}:={}^o\!(\gamma_{f,g})^{-1}
\qquad\text{and}\qquad
{}^o\!\delta_{{}^o\!A}:={}^o\!(\delta_A)^{-1}
$$
whose composition and unit axioms are again derived
straightforwardly from the same for $F$.

(v)\ \
Let $\cC$ be a third $2$-category, and $G:\cB\to\cC$ another
pseudo-functor. We may define a composition $G\circ F:\cA\to\cB$,
which is the pseudo-functor such that
$$
G\circ F(A):=G(FA)
\qquad\text{and}\qquad
(G\circ F)_{AB}:=G_{FA,FB}\circ F_{AB}
\qquad
\text{for every $A,B\in\Ob(\cA)$}.
$$
Denote by $(\delta^F,\gamma^F)$ and $(\delta^G,\gamma^G)$
the coherence constraints of $F$ and respectively $G$; the
coherence constraint of $G\circ F$ is then the pair
$(\delta^{G\circ F},\gamma^{G\circ F})$ such that
$$
\delta^{G\circ F}_A:=G(\delta^F_A)\odot\delta_{FA}^G
\qquad
\gamma^{G\circ F}_{f,g}:=G(\gamma^F_{f,g})\odot\gamma^G_{Ff,Fg}
$$
for every $A,B,C\in\Ob(\cA)$ and every $1$-cell $f$ of
$\cA(A,B)$ and $g$ of $\cA(B,C)$. Let us verify the
composition axiom for $\gamma^{G\circ F}$ : given
a composable sequence of $1$-cells $f,g,h$, we have
$$
\begin{aligned}
\gamma^{G\circ F}_{g\circ f,h}\odot(GFh*\gamma^{G\circ F}_{f,g})
=\, & G(\gamma^F_{g\circ f,h})\odot\gamma^G_{F(g\circ f),Fh}\odot
(GFh*(G(\gamma^F_{f,g})\odot\gamma^G_{Ff,Fg})) \\
=\, & G(\gamma^F_{g\circ f,h})\odot\gamma^G_{F(g\circ f),Fh}\odot
(GFh*G(\gamma^F_{f,g}))\odot(GFh*\gamma^G_{Ff,Fg}) \\
=\, & G(\gamma^F_{g\circ f,h})\odot G(Fh*\gamma^F_{f,g})\odot
\gamma^G_{Fg\circ Ff,Fh}\odot(GFh*\gamma^G_{Ff,Fg}) \\
=\, & G(\gamma^F_{g\circ f,h}\odot(Fh*\gamma^F_{f,g}))\odot
\gamma^G_{Ff,Fh\circ Fg}\odot(\gamma^G_{Fg,Fh}*GFf) \\
=\, & G(\gamma^F_{f,h\circ g}\odot(\gamma^F_{g,h}*Ff))\odot
\gamma^G_{Ff,Fh\circ Fg}\odot(\gamma^G_{Fg,Fh}*GFf) \\
=\, & G(\gamma^F_{f,h\circ g})\odot G(\gamma^F_{g,h}*Ff)\odot
\gamma^G_{Ff,Fh\circ Fg}\odot(\gamma^G_{Fg,Fh}*GFf) \\
=\, & G(\gamma^F_{f,h\circ g})\odot\gamma^G_{Ff,F(h\circ g)}\odot
(G(\gamma^F_{g,h})*GFf)\odot(\gamma^G_{Fg,Fh}*GFf) \\
=\, & \gamma^{G\circ F}_{f,h\circ g}\odot
((G(\gamma^F_{g,h})\odot\gamma^G_{Fg,Fh})*GFf) \\
=\, & \gamma^{G\circ F}_{f,h\circ g}\odot(\gamma^{G\circ F}_{g,h}*GFf)
\end{aligned}
$$
where the third and seventh equalities follow from (ii) (applied
to $G$), the fourth from (i) and the composition axiom for $G$,
the fifth from the composition axiom for $F$, the sixth from (i).

Next, to check the unit axiom, we compute :
$$
\begin{aligned}
\gamma_{\one_A,f}^{G\circ F}\odot(GFf*\delta^{G\circ F}_A)
=\, & G(\gamma^F_{\one_A,f})\odot\gamma^G_{F\one_A,Ff}\odot
(GFf*(G(\delta^F_A)\odot\delta^G_{FA})) \\
=\, & G(\gamma^F_{\one_A,f})\odot\gamma^G_{F\one_A,Ff}\odot
(GFf*G(\delta^F_A))\odot(GFf*\delta^G_{FA}) \\
=\, & G(\gamma^F_{\one_A,f})\odot G(Ff*\delta^F_A)\odot
\gamma^G_{\one_{FA},Ff}\odot(GFf*\delta^G_{FA}) \\
=\, & G(\gamma^F_{\one_A,f}\odot(Ff*\delta^F_A)) \\
=\, & G(\one_{Ff}) \\
=\, & \one_{GFf}
\end{aligned}
$$
where the third equality follows from (ii), the fourth from
(i) and the unit axiom for $G$, and the fifth from the unit
axiom for $F$. Similarly one verifies the commutativity of
the remaining diagram required for the unit axiom : details
left to the reader.

(vi)\ \
Also, if $H:\cC\to\cD$ is any other pseudo-functor, then
$H\circ(G\circ F)=(H\circ G)\circ F$. Indeed, taking into
account (i) we may compute :
$$
\begin{aligned}
\gamma^{H\circ(G\circ F)}_{f,g}=\, &
H(\gamma^{G\circ F}_{f,g})\odot\gamma^H_{G\circ Ff,G\circ Fg} \\
=\, & H\circ G(\gamma^F_{f,g})\odot H(\gamma^G_{Ff,Gf})\odot
\gamma^H_{G\circ Ff,G\circ Fg} \\
=\, & H\circ G(\gamma^F_{f,g})\odot\gamma^{H\circ G}_{Ff,Gf} \\
=\, & \gamma^{(H\circ G)\circ F}_{f,g}.
\end{aligned}
$$
Similarly one checks that
$\delta^{H\circ(G\circ F)}_A=\delta^{(H\circ G)\circ F}_A$ for
every $A\in\Ob(\cA)$, whence the assertion.

(vii)\ \
Notice as well that if $F$ and $G$ are strict, then the same
holds for $G\circ F$.
\end{remark}

\begin{remark}\label{rem_remark-cat-cats}
(i)\ \
Remark \ref{rem_pseudo-funct}(v) shows that the set
of all $\sU$-small $2$-categories together with all
pseudo-functors between them forms a category with
small $\Hom$-sets that we denote
$$
\sU\tdu 2\tdu\bCat
$$
(and as always, we shall usually drop the universe $\sU$
from the notation). We do not know whether $2\tdu\bCat$
is a complete category : for instance, if $F:\cA\to\cB$
and $F':\cA'\to\cB$ are any two pseudo-functors, it is
not clear whether the fibre product of $F$ and $F'$ is
always representable in $2\tdu\bCat$. However, in case
$F$ and $F'$ are both strict, this fibre product is
indeed representable by a $2$-category that we denote
$$
\cA\times_{(F,F')}\cA'
$$
(or simply $\cA\times_\cB\cA'$, if the notation is not
ambiguous). Namely, the objects of $\cA\times_\cB\cA'$
are all the pairs $(A,A')$ consisting of objects $A$ of
$\cA$ and $A'$ of $\cA'$ such that $FA=F'A'$, and the
$\Hom$-categories are given by the rule :
$$
\cA\times_\cB\cA'((A_1,A'_1),(A_2,A'_2)):=
\cA(A_1,A_2)\times_{\cB(FA_1,FA_2)}\cA(A_2,A'_2)
$$
for every pair of objects $(A_1,A'_1),(A_2,A'_2)$ of
$\cA\times_\cB\cA'$ (see example \ref{ex_cat-cats}(i)).
The composition laws and unit functors for
$\cA\times_\cB\cA'$ are deduced from those of $\cA$ and
$\cA'$, in the obvious fashion. A universal cone for
this fibre product is provided by the two projections
$\cA\leftarrow\cA\times_\cB\cA'\to\cA'$, {\em i.e.} the
strict pseudo-functors defined in the obvious fashion.

(ii)\ \
For instance, just as in \eqref{subsec_Morph-cat}, for every
$2$-category $\cA$ we have two natural {\em source}
and {\em target} strict pseudo-functors :
$$
\cA\xleftarrow{\ \ss\ }2\tdu\sMorph(\cA)\xrightarrow{\ \st\ }\cA
$$
(notation of example \ref{ex_2-arrows}(i)). Namely, $\ss$
assigns to every $1$-cell $f:A\to B$ the object $A$ of
$\cA$, to every diagram \eqref{eq_2-cell} the $1$-cell
$g:\ss(f)\to\ss(f')$ of $\cA$, and to every $2$-cell
$(\beta,\gamma)$ of $2\tdu\sMorph(\cA)$, the $2$-cell
$\beta$ of $\cA$. The reader may likewise spell out
the definition of $\st$. Then we have as well a natural
{\em composition} strict pseudo-functor
$$
\sc:2\tdu\sMorph(\cA)\times_{(\st,\ss)}2\tdu\sMorph(\cA)
\to 2\tdu\sMorph(\cA)
$$
which assigns :
\begin{itemize}
\item
to every pair of $1$-cells $(f,g)$ such
that $\st(f)=\ss(g)$, the composition $\sc(f,g):=g\circ f$
\item
to every pair of $1$-cells of $2\tdu\sMorph(\cA)$ :
\set\begin{equation}\label{eq_two-two-cells}
{\diagram A \ar[r]^-f \ar[d]_{h}
\drtwocell\omit{_\alpha\ } & B \ar[d]|k \ar[r]^-g
\drtwocell\omit{_\beta\ } & C \ar[d]^i \\
A' \ar[r]_-{f'} & B' \ar[r]_-{g'} & C'
\enddiagram}
\end{equation}
the $1$-cell of $2\tdu\sMorph(\cA)$
$$
\sc((\alpha,h,k),(\beta,k,i)):=(\beta\boxminus\alpha,h,i)
$$
\item
for $(f,g),(f',g')$ as in \eqref{eq_two-two-cells}, and
every $2$-cells
$(\phi_1,\phi_2):(\alpha_1,h_1,k_1)\Rightarrow(\alpha_2,h_2,k_2)$
in $2\tdu\sMorph(f,f')$ and
$(\phi_2,\phi_3):(\beta_1,k_1,i_1)\Rightarrow(\beta_2,k_2,i_2)$
in $2\tdu\sMorph(g,g')$, the pair
$$
\sc((\phi_1,\phi_2),(\phi_2,\phi_3)):=(\phi_1,\phi_3)
$$
which is a well defined $2$-cell
$\sc((\alpha_1,h_1,k_1),(\beta_1,k_1,i_1))\Rightarrow
\sc((\alpha_2,h_2,k_2),(\beta_2,k_2,i_2))$ in
$2\tdu\sMorph(g\circ f,g'\circ f')$, since by assumption
we have
$$
\frac{\alpha_1}{h_1|g_1}\boxminus\frac{\phi_1}{g_1|g_2}=
\frac{\phi_2}{h_1|h_2}\boxminus\frac{\alpha_2}{h_2|g_2}
\qquad
\frac{\beta_1}{k_1|h_1}\boxminus\frac{\phi_2}{h_1|h_2}=
\frac{\phi_3}{k_1|k_2}\boxminus\frac{\beta_2}{k_2|h_2}
$$
and therefore
$$
\begin{aligned}
(\beta_1\boxminus\alpha_1)\boxminus\phi_1=\, &
\beta_1\boxminus(\alpha_1\boxminus\phi_1) \\
=\, & \beta_1\boxminus(\phi_2\boxminus\alpha_2) \\
=\, & (\beta_1\boxminus\phi_2)\boxminus\alpha_2 \\
=\, & (\phi_3\boxminus\beta_2)\boxminus\alpha_2 \\
=\, & \phi_3\boxminus(\beta_2\boxminus\alpha_2).
\end{aligned}
$$
\end{itemize}
With these rules, the functoriality of $\sc$ for
$2$-cells is immediate. To check the functoriality
of $\sc$ for $1$-cells, consider a diagram as in
proposition \ref{prop_square-algebra}(ii), and to
ease notation, set $\underline\alpha:=(\alpha,h,k)$
and define likewise $\underline\alpha'$, $\underline\beta$
and $\underline\beta'$. Then proposition
\ref{prop_square-algebra}(ii) translates as the identity
$$
\sc(\underline\alpha'\circ\underline\alpha,
\underline\beta'\circ\underline\beta)=
\sc(\underline\alpha',\underline\beta')\circ
\sc(\underline\alpha,\underline\beta)
$$
as required.
\end{remark}

\begin{example}\label{ex_2-funct-of-arrows}
(i)\ \
Let $\cA$, $\cB$ be two $2$-categories, and $F:\cA\to\cB$
any pseudo-functor, with coherence constraint
$(\gamma^F,\delta^F)$. Then $F$ induces a pseudo-functor
$$
2\tdu\sMorph(F):2\tdu\sMorph(\cA)\to 2\tdu\sMorph(\cB)
$$
by the following rule. To every $1$-cell $f:A\to B$ in
$\cA$ we assign the $1$-cell $Ff:FA\to FB$. To every
diagram \eqref{eq_2-cell} we assign the diagram
$$
{\diagram FA \ar[r]^-{Ff} \ar[d]_{Fg}
\drtwocell\omit{_\ \alpha^F} & FB \ar[d]^{Fh} \\
FA' \ar[r]_-{Ff'} & FB'
\enddiagram}
\qquad
\text{where
$\alpha^F:=(\gamma^F_{g,f'})^{-1}\odot F\alpha\odot\gamma^F_{f,h}$}
$$
and to any $2$-cell
$(\beta,\delta):(g_1,h_1,\alpha_1)\Rightarrow(g_2,h_2,\alpha_2)$
in $2\tdu\sMorph(f,f')$, we assign the pair $(F\beta,F\delta)$.
Indeed, let us check that the latter is a well defined
$2$-cell in $2\tdu\sMorph(Ff,Ff')$; the assertion boils
down to the identity :
$$
(\gamma^F_{g_2,f'})^{-1}\odot
F\alpha_2\odot\gamma^F_{f,h_2}\odot(F\delta*Ff)=
(Ff'*F\beta)\odot(\gamma_{g_1,f'})^{-1}\odot
F\alpha_1\odot\gamma^F_{f,h_1}.
$$
However, from remark \ref{rem_pseudo-funct}(ii) we get
$$
(Ff'*F\beta)\odot(\gamma^F_{g_1,f'})^{-1}=
(\gamma^F_{g_2,f'})^{-1}\odot F(f'*\beta)
$$
hence it suffices to check that
$$
F\alpha_2\odot\gamma^F_{f,h_2}\odot(F\delta*Ff)=
F(f'*\beta)\odot F\alpha_1\odot\gamma^F_{f,h_1}.
$$
But by assumption, we know that
$(f'*\beta)\odot\alpha_1=\alpha_2\odot(\delta*f)$, so
-- taking into account remark \ref{rem_pseudo-funct}(i) --
we are further reduced to showing that
$$
\gamma^F_{f,h_2}\odot(F\delta*Ff)=F(\delta*f)\odot\gamma^F_{f,h_1}
$$
which follows again from remark \ref{rem_pseudo-funct}(ii).
Lastly, the coherence constraint of $2\tdu\sMorph(F)$ assigns
to any diagram \eqref{eq_unit-in-morph} (resp.
\eqref{eq_star-comp}) the pair
\set\begin{equation}\label{eq_double-constraint}
(\delta^F_A,\delta^F_B)
\qquad
\text{(resp.\ \ $(\gamma^F_{g,k},\gamma^F_{h,i})$\ )}.
\end{equation}
Indeed, let us check that $(\delta^F_A,\delta^F_B)$ is
a well defined $2$-cell in $2\tdu\sMorph(Ff,Ff)$ : the
assertion boils down to the identity
$$
(Ff*\delta^F_A)\odot\one^F_f=\one^F_f\odot(\delta^F_B*Ff)
\qquad\text{where}\qquad
\one^F_f:=(\gamma^F_{\one_A,f})^{-1}\odot\gamma^F_{f,\one_B}
$$
which in turns follows easily from the unit axiom for
$\delta^F$. Likewise, the assertion that
$(\gamma^F_{g,k},\gamma^F_{h,i})$ is a well defined $2$-cell
in $2\tdu\sMorph(Ff,Ff'')$ comes down to the identity
\set\begin{equation}\label{eq_damn-overfull}
\gamma^F_{g,k}\boxminus(\beta^F\boxvert\alpha^F)=
(\beta\boxvert\alpha)^F\boxminus\gamma^F_{h,i}.
\end{equation}
However, by definition, the left-hand side of
\eqref{eq_damn-overfull} equals
$$
(Ff''*\gamma^F_{g,k})\odot((\gamma^F_{k,f''})^{-1}*Fg)
\odot(F\beta*Fg)\odot(\gamma^F_{f',i}*Fg)\odot
(Fi*(\gamma^F_{g,f'})^{-1})\odot(Fi*(F\alpha\odot\gamma^F_{f,h}))
$$
and the composition axiom for $F$ yields the identities :
$$
\begin{aligned}
(Ff''*\gamma^F_{g,k})\odot((\gamma^F_{k,f''})^{-1}*Fg)=\, &
(\gamma^F_{k\circ g,f''})^{-1}\odot\gamma^F_{g,f''\circ k} \\
(\gamma^F_{f',i}*Fg)\odot(Fi*(\gamma^F_{g,f'})^{-1})=\, &
(\gamma^F_{g,i\circ f'})^{-1}\odot\gamma^F_{f'\circ g,i}.
\end{aligned}
$$
So the left-hand side of \eqref{eq_damn-overfull} also equals
$$
(\gamma^F_{k\circ g,f''})^{-1}\odot\gamma^F_{g,f''\circ k}\odot
(F\beta*Fg)\odot(\gamma^F_{g,i\circ f'})^{-1}\odot
\gamma^F_{f'\circ g,i}\odot(Fi*(F\alpha\odot\gamma^F_{f,h})).
$$
Next, remark \ref{rem_pseudo-funct}(ii) yields the identities
$$
\begin{aligned}
(F\beta*Fg)\odot(\gamma^F_{g,i\circ f'})^{-1}=\, &
(\gamma^F_{g,f''\circ k})^{-1}\odot F(\beta*g) \\
\gamma^F_{f'\circ g,i}\odot(Fi*F\alpha)=\, &
F(i*\alpha)\odot\gamma^F_{h\circ f,i}.
\end{aligned}
$$
So we have to show that the right-hand side of
\eqref{eq_damn-overfull} equals
$$
(\gamma^F_{k\circ g,f''})^{-1}\odot F(\beta*g)\odot
F(i*\alpha)\odot\gamma^F_{h\circ f,i}\odot(Fi*\gamma^F_{f,h}).
$$
But by unwinding the definition, we see that the
right-hand side of \eqref{eq_damn-overfull} equals
$$
(\gamma^F_{k\circ g,f''})^{-1}\odot F(\beta\boxvert\alpha)
\odot\gamma^F_{f,i\circ h}\odot(\gamma^F_{h,i}*Ff)=
(\gamma^F_{k\circ g,f''})^{-1}\odot F(\beta\boxvert\alpha)
\odot\gamma^F_{h\circ f,i}\odot(Fi*\gamma^F_{f,h})
$$
(where the equality holds by the composition axiom). Thus,
we come down to checking that
$$
F(\beta*g)\odot F(i*\alpha)=
F(\beta\boxvert\alpha)
$$
which follows from remark \ref{rem_pseudo-funct}(i).
With these rules, the functoriality of $2\tdu\sMorph(F)$
is clear from remark \ref{rem_pseudo-funct}(i), and the
naturality of \eqref{eq_double-constraint}, as well as
the composition and unit axioms for the latter, follow
from the respective properties of $(\gamma^F,\delta^F)$,
by a direct inspection.

(ii)\ \
In the situation of (i), let $f:A\to B$ be any $1$-cell
of $\cA$. Then the unit axiom for the coherence constraint
$(\gamma^F,\delta^F)$ implies the identity :
$$
\one_f^F=(Ff*\delta^F_A)\odot((\delta^F_B)^{-1}*Ff).
$$

(iii)\ \
Let $\cC$ be a third $2$-category, and $G:\cB\to\cG$
another pseudo-functor; then we have
$$
2\tdu\sMorph(G\circ F)=2\tdu\sMorph(G)\circ 2\tdu\sMorph(F)
$$
(details left to the reader).
\end{example}

\begin{example}\label{ex_change-orientation}
Let $\cA$ be any $2$-category.

(i)\ \
We have a natural strict isomorphism of $2$-categories
$$
{}^o(2\tdu\sMorph(\cA))^o\isom
2\tdu\sMorph({}^o\cA^o)
$$
that assigns :
\begin{itemize}
\item
to every object $f:A\to B$ of $2\tdu\sMorph(\cA)$
the object ${}^of^o:{}^oB^o\to{}^oA^o$
\item
to every oriented square \eqref{eq_2-cell} the corresponding
oriented square ${}^o\eqref{eq_2-cell}^o$, which is a $1$-cell
${}^o\alpha^o:{}^of'^o\to{}^of^o$ in $2\tdu\sMorph({}^o\cA^o)$
\item
to every $2$-cell $(\beta,\gamma)$ in $2\tdu\sMorph(\cA)(f,f')$
the pair $({}^o\gamma^o,{}^o\beta^o)$ which is a $2$-cell in
$2\tdu\sMorph({}^o\cA^o)({}^of'^o,{}^of^o)$.
\end{itemize}

(ii)\ \
We deduce from (i) a strict isomorphism of $2$-categories
$$
(2\tdu\sMorph(\cA^o))^o\isom {}^o(2\tdu\sMorph({}^o\cA)).
$$
These latter $2$-categories can be described as follows.
The objects are the same as those of $2\tdu\sMorph(\cA)$,
and the $1$-cells $f\to f'$ are the oriented squares just
like \eqref{eq_2-cell}, but {\em with reversed direction
of the orienting arrow}, {\em i.e.} $\alpha$ is replaced
by a $2$-cell $f'\circ g\Rightarrow h\circ f$. Given
two $1$-cells $(g_1,h_1,\alpha_1),(g_2,h_2,\alpha_2):f\to f'$
in $2\tdu\sMorph(\cA^o)^o(f,f')$, the $2$-cells
$(g_1,h_1,\alpha_1)\Rightarrow(g_2,h_2,\alpha_2)$ are
still all the pairs
$(\beta:g_1\Rightarrow g_2,\gamma:h_1\Rightarrow h_é)$
such that $\alpha_1\boxminus\beta=\gamma\boxminus\alpha_2$,
and the composition rule for such $2$-cells is the same
as in $2\tdu\sMorph(\cA)$. However, the composition law
for $1$-cells in $2\tdu\sMorph(\cA^o)^o$ is given by the
join operation $\boxminus$ instead of $\boxvert$ : details
left to the reader. 
\end{example}

\begin{example}\label{ex_opposite-cat-as-pseudo-fun}
We have a natural strict isomorphism of $2$-categories :
$$
(-)^o:{}^o\bCat\isom\bCat
$$
that assigns to every small category $\cC$ the opposite
category $\cC^o$, to every functor $F:\cB\to\cC$ the
opposite functor $F^o:\cB^o\to\cC^o$, and to every natural
transformation $\beta:F\Rightarrow G$ the opposite
transformation $\beta^o:G^o\Rightarrow F^o$.
\end{example}

\subsection{Pseudo-natural transformations and their modifications}
The following definition introduces the $2$-categorical analogue
of a natural transformation between functors.

\begin{definition}\label{def_pseudo-natural}
Consider two pseudo-functors $F,G:\cA\to\cB$ between
$2$-categories $\cA$, $\cB$.

(i)\ \
A {\em lax-natural transformation\/} $\alpha:F\Rightarrow G$
is the datum of :

$\bullet$\ \
For every object $A$ of $\cA$, a $1$-cell $\alpha_A:FA\to GA$.

$\bullet$\ \
For every $A,B\in\Ob(\cA)$, a natural transformation
$\tau_{AB}$ between two functors $\cA(A,B)\to\cB(FA,GB)$, as
shown by the (not necessarily commutative) diagram
\set\begin{equation}\label{eq_pseudo-nat}
{\diagram \cA(A,B) \ar[d]_{G_{AB}} \ar[rrr]^-{F_{AB}} &
\drtwocell\omit{^\tau_{AB}\ \ \ } & &
\cB(FA,FB) \ar[d]^{H_\cB(\one_{FA},\alpha_B)} \\
\cB(GA,GB) \ar[rrr]_-{H_\cB(\alpha_A,\one_{GB})} & & &
\cB(FA,GB)
\enddiagram}\end{equation}
where $H_\cB(\one_{FA},\alpha_B)$ and $H_\cB(\alpha_A,\one_{GB})$
are as in example \ref{ex_first-from-Cats-to-PsFuns}.
The datum $\tau_{\bullet\bullet}$ is called the
{\em coherence constraint\/} for $\alpha$, and is required
to satisfy the following {\em coherence axioms} (in which
we denote by $(\delta^F,\gamma^F)$ and $(\delta^G,\gamma^G)$
the coherence constraints for $F$ and respectively $G$) :
\begin{itemize}
\item
For every $A\in\Ob(\cA)$, the following diagram commutes :
$$
\xymatrix{
\alpha_A \ar@{=>}[r]^-{\one_{\alpha_{\!A}}}
\ar@{=>}[d]_{\one_{\alpha_{\!A}}} &
\one_{GA}\circ\alpha_A
\ar@{=>}[rr]^-{\delta^G_{\!A}*\alpha_{\!A}} & &
G(\one_A)\circ\alpha_A \ar@{=>}[d]^{\tau_{\one_{\!A}}} \\
\alpha_A\circ\one_{FA}
\ar@{=>}[rrr]^-{\alpha_{\!A}*\delta^F_{\!A}} & & &
\alpha_A\circ F(\one_A).
}$$
\item
For each pair of $1$-cells $A\xrightarrow{f}B\xrightarrow{g}C$
in $\cA$, the following diagram commutes :
$$
\xymatrix{
Gg\circ Gf\circ\alpha_A \ar@{=>}[rr]^-{Gg*\tau_f}
\ar@{=>}[d]_{\gamma^G_{f,g}*\alpha_{\!A}} & &
Gg\circ\alpha_B\circ Ff \ar@{=>}[rr]^-{\tau_g*Ff} & &
\alpha_C\circ Fg\circ Ff \ar@{=>}[d]^{\alpha_C*\gamma^F_{f,g}} \\
G(g\circ f)\circ\alpha_A \ar@{=>}[rrrr]^-{\tau_{g\circ f}} & & & &
\alpha_C\circ F(g\circ f).
}$$
\end{itemize}

(ii)\ \
A lax-natural transformation $\alpha$ as in (i) is called
a {\em pseudo-natural transformation} (resp. {\em strict})
if $\tau_{AB}$ is an isomorphism of functors (resp. if
\eqref{eq_pseudo-nat} commutes and $\tau_{AB}$ is the
identity natural transformation) for every $A,B\in\Ob(\cA)$.

(iii)\ \
A pseudo-natural transformation $\alpha$ as in (i) is
called a {\em pseudo-natural isomorphism} if $\alpha_A$
is an invertible $1$-cell, for every $A\in\Ob(\cA)$.
\end{definition}

\begin{example}\label{ex_4-funct-of-arrows}
For any $2$-category $\cC$, denote by
$$
p_\cC,p'_\cC:2\tdu\sMorph(\cC)\times_{(\st,\ss)}
2\tdu\sMorph(\cC)\to 2\tdu\sMorph(\cC)
$$
the projections, {\em i.e.} the strict pseudo-functors that
yield a universal cone for this fibre product, where $\st$
and $\ss$ are the target and source pseudo-functors for
$2\tdu\sMorph(\cC)$ as in remark \ref{rem_remark-cat-cats}(ii).
Now, let $\cA,\cB$ be two $2$-categories, and $F:\cA\to\cB$ a
pseudo-functor; it is easily seen that
$$
\st\circ 2\tdu\sMorph(F)\circ p_\cA=
\ss\circ 2\tdu\sMorph(F)\circ p'_\cA
$$
where $2\tdu\sMorph(F)$ is defined as in example
\ref{ex_2-funct-of-arrows}(i). There follows a unique
pseudo-functor
$$
2\tdu\sMorph(F,F):
2\tdu\sMorph(\cA)\times_{(\st,\ss)}
2\tdu\sMorph(\cA)\to
2\tdu\sMorph(\cB)\times_{(\st,\ss)}
2\tdu\sMorph(\cB)
$$
such that
$$
p_\cB\circ 2\tdu\sMorph(F,F)=2\tdu\sMorph(F)\circ p_\cA
\qquad\text{and}\qquad
p'_\cB\circ 2\tdu\sMorph(F,F)=2\tdu\sMorph(F)\circ p'_\cA
$$
and we claim that the coherence constraint $\gamma^F$ of
$F$ can be viewed as a strict pseudo-natural isomorphism
$$
\xymatrix{
2\tdu\sMorph(\cA)\times_{(\st,\ss)}2\tdu\sMorph(\cA)
\ar[d]_{2\tdu\sMorph(F,F)} \ar[rrr]^-\sc & 
\drtwocell\omit{^\gamma^F\ \ } & & 2\tdu\sMorph(\cA)
\ar[d]^{2\tdu\sMorph(F)} \\
2\tdu\sMorph(\cB)\times_{(\st,\ss)}2\tdu\sMorph(\cB)
\ar[rrr]_-\sc & & & 2\tdu\sMorph(\cB).
}$$
Namely, for every pair of composable $1$-cells $(f,g)$
of $\cA$, we regard
$\gamma^F_{f,g}:Fg\circ Ff\Rightarrow F(g\circ f)$ as
a $1$-cell of $2\tdu\sMorph(\cB)$. Then the strictness
of $\gamma^F$ boils down to the identity
$$
\gamma^F_{f',g'}\boxvert(\beta^F\boxminus\alpha^F)=
(\beta\boxminus\alpha)^F\boxvert\gamma^F_{f,g}
$$
for every diagram \eqref{eq_two-two-cells} in $\cA$.
For the verification, notice that the associativity
constraint for the functor $(2\tdu\sMorph(F^o))^o$
assigns to \eqref{eq_two-two-cells} the pair
$(\gamma_{f,g}^o,\gamma_{f',g'}^o)^o=(\gamma_{f,g},\gamma_{f',g'})$
(see remark \ref{rem_pseudo-funct}(iv) and example
\ref{ex_2-funct-of-arrows}(i)), and then the foregoing
identity just translates the assertion that
$(\gamma_{f,g},\gamma_{f',g'})$ is a well defined $2$-cell
in $2\tdu\sMorph(F^oh^o,F^oi^o)^o$ (see example
\ref{ex_change-orientation}(ii)).
The coherence axioms for $\gamma^F$ follow easily,
by unwinding the definitions.
\end{example}

\begin{example}\label{ex_2-cats}
(i)\ \ 
Any category $\cA$ can be regarded as a $2$-category in a
natural way : namely, for any two objects $A$ and $B$ of $\cA$
one lets $\cA(A,B)$ be the discrete category $\Hom_\cA(A,B)$;
hence the only $2$-cells of $\cA$ are the identities
$\one_f:f\Rightarrow f$, for every morphism $f:A\to B$. The
composition bifunctor $c_{ABC}$ is of course given (on $1$-cells)
by the composition law for morphisms of $\cA$. Likewise, the
functor $u_A$ assigns to every object $A$ its identity endomorphism.

(ii)\ \ 
In the same vein, every functor between usual categories
is a strict pseudo-functor between the corresponding
$2$-categories as in (i). Finally, every natural transformation
of usual functors can be regarded naturally as a strict
pseudo-natural transformation between the corresponding
strict pseudo-functors.
\end{example}

\begin{remark}\label{rem_pseudo-natural}
With the notation of definition \ref{def_pseudo-natural},
we have :

(i)\ \
The coherence constraint of $\alpha$ can be regarded
as a system of oriented squares
$$
{\diagram FA \ar[r]^-{\alpha_A} \ar[d]_{Ff}
\drtwocell\omit{_\ \ \tau_f} & GA \ar[d]^{Gf} \\
FB \ar[r]_-{\alpha_B} & GB
\enddiagram}
\qquad
\text{for every $1$-cell $f:A\to B$ of $\cA$}.
$$
Then $\alpha$ is a pseudo-natural transformation if and only
if the orientation $\tau_f$ is invertible for every $1$-cell
$f$ of $\cA$. The naturality of $\tau_{\bullet\bullet}$ comes
down to the commutativity of the diagram
$$
\xymatrix{ Gf\circ\alpha_A \ar@{=>}[r]^-{\tau_f}
\ar@{=>}[d]_{G\beta*\alpha_A} &
\alpha_B\circ Ff \ar@{=>}[d]^{\alpha_B*F\beta} \\
Gf'\circ\alpha_A \ar@{=>}[r]^-{\tau_{f'}} & \alpha_B\circ Ff'
}$$
for every $A,B\in\Ob(\cA)$ and every $2$-cell
$\beta:f\Rightarrow f'$ in $\cA(A,B)$. The latter
can be also be written as the identity
$$
\tau_f\boxminus F\beta=G\beta\boxminus\tau_{f'}.
$$
Likewise, the coherence axioms can be written as the
identities
$$
\delta^G_A\boxminus\tau_{\one_A}=\one_{\alpha_A}\boxminus\delta^F_A
\qquad
(\tau_g\boxvert\tau_f)\boxminus\gamma^F_{f,g}=
\gamma^G_{f,g}\boxminus\tau_{g\circ f}
$$
for every $A\in\Ob(\cA)$ and every composable pair
$(f,g)$ of $1$-cells of $\cA$.

(ii)\ \
The foregoing list of requirements becomes more intelligible,
when we observe that a lax-natural transformation
$\alpha:F\Rightarrow G$ is equivalent to the datum of
a pseudo-functor
$$
\tilde\alpha:\cA\to 2\tdu\sMorph(\cB)
\qquad\text{such that}\qquad
\ss\circ\tilde\alpha=F
\qquad\text{and}\qquad
\st\circ\tilde\alpha=G
$$
and $\alpha$ is a pseudo-natural transformation if and only
if $\tilde\alpha(f)$ is an invertible $2$-cell for every
$1$-cell $f$ of $\cA$. Namely, in light of (i) we obtain such
$\tilde\alpha$ by setting
$$
\tilde\alpha(A):=\alpha_A
\qquad
\tilde\alpha(f):=\tau_f
\qquad
\tilde\alpha(\beta):=(F\beta,G\beta)
$$
for every $A\in\Ob(\cA)$, every $1$-cell $f$ of $\cA$,
and every $2$-cell $\beta$ of $\cA$. The coherence
constraint $(\delta,\gamma)$ of $\tilde\alpha$ are
defined by setting
$$
\delta_A:=(\delta^F_A,\delta^G_A)
\qquad
\gamma_{f,g}:=(\gamma^F_{f,g},\gamma^G_{f,g})
$$
for every $A\in\Ob(\cA)$ and every composable
pair $(f,g)$ of $1$-cells of $\cA$. Especially,
$\tilde\alpha$ is a strict pseudo-functor if and
only if the same holds for both $F$ and $G$.

(iii)\ \
If $\alpha:F\Rightarrow G$ and $\beta:G\Rightarrow H$
are two lax-natural transformations, let
$\tilde\alpha,\tilde\beta:\cA\to 2\tdu\sMorph(\cB)$ be
the associated pseudo-functors, as in (ii); then we may
proceed as in example \ref{ex_cat-cats}(iii) to define
the composition
$$
\beta\odot\alpha:F\Rightarrow H.
$$
Namely, it shall be the lax-natural transformation
associated to the pseudo-functor
$$
\cA\xrightarrow{(\tilde\alpha,\tilde\beta)}
2\tdu\sMorph(\cA)\times_{(\st,\ss)}2\tdu\sMorph(\cA)
\xrightarrow{\ \sc\ } 2\tdu\sMorph(\cA)
$$
(notation of remark \ref{rem_remark-cat-cats}(ii)).
Unwinding the definitions, we see that this is the
lax-natural transformation given by the rule
$$
A\mapsto\beta_A\circ\alpha_A
\qquad
\text{for every $A\in\Ob(\cA)$}
$$
and with coherence constraint
$\tau^{\beta\odot\alpha}_{\bullet\bullet}$ given by the rule
(notation of remark \ref{rem_remark-cat-cats}(ii)) :
$$
(f:A\to B)\mapsto
\tau^\beta_f\boxminus\tau^\alpha_f
\qquad
\text{for every $1$-cell $f$ of $\cA$}
$$
where $\tau^\alpha_{\bullet\bullet}$ (resp. $\tau^\beta_{\bullet\bullet}$)
denotes the coherence constraint of $\alpha$ (resp. of $\beta$).
Hence, if both $\alpha$ and $\beta$ are pseudo-natural (resp.
strict), the same holds for $\beta\odot\alpha$. Since the
composition law $\boxminus$ is associative (proposition
\ref{prop_square-algebra}(i) and remark \ref{rem_square-algebra}),
we deduce that
$$
\lambda\odot(\beta\odot\alpha)=(\lambda\odot\beta)\odot\alpha
$$
for every lax-natural transformation
$\lambda:H\Rightarrow K$ (and any pseudo-functor
$K:\cA\to\cB$).

(iv)\ \
If $\cC$ is a third $2$-category, $F',G':\cB\to\cC$
two other pseudo-functors, and $\alpha':F'\Rightarrow G'$
another lax-natural transformation, we may likewise
define the {\em Godement products}
$$
\alpha'*F:F'\circ F\Rightarrow G'\circ F
\qquad
G'*\alpha:G'\circ F\Rightarrow G'\circ G.
$$
Namely, they shall be the pseudo-natural transformations
attached to the pseudo-functors
$$
\tilde\alpha'\circ F
\qquad\text{and respectively}\qquad
2\tdu\sMorph(G')\circ\tilde\alpha
$$
where $\tilde\alpha':\cB\to 2\tdu\sMorph(\cC)$ is the
pseudo-functor associated with $\alpha'$, as in (i),
and we define $2\tdu\sMorph(G')$ as in example
\ref{ex_2-funct-of-arrows}(i). If $\alpha$
(resp. $\alpha'$) is pseudo-natural, the same holds
for $G'*\alpha$ (resp. for $\alpha'*F$). In view of
example \ref{ex_2-funct-of-arrows}(iii) and remark
\ref{rem_pseudo-funct}(vi), it is then clear that
$$
\alpha'*(F\circ K)=(\alpha'*F)*K
\qquad\text{and}\qquad
(K'\circ G')*\alpha=K'*(G'*\alpha)
$$
for any other $2$-category $\cD$, and every pseudo-functors
$K:\cD\to\cA$ and $K':\cC\to\cD$. Furthermore, if $\beta$
is as in (ii), and $\beta':G'\Rightarrow H'$ is another
lax-natural transformation (to a target pseudo-functor
$H':\cB\to\cC$), then we get straightforwardly :
$$
(G'*\alpha)*K=G'*(\alpha*K)
\qquad\text{and}\qquad
(\beta'\odot\alpha')*F=(\beta'*F)\odot(\alpha'*F).
$$
On the other hand, example \ref{ex_4-funct-of-arrows}
says that the coherence constraint $\gamma^{G'}$ of $G'$
induces a pseudo-natural isomorphism from the pseudo-functor
associated to $(G'*\beta)\odot(G'*\alpha)$ to the
pseudo-functor associated to $G'*(\beta\odot\alpha)$.
Thus, if $G'$ is a strict pseudo-functor, then these
two lax-natural transformations coincide, but in
general they will be different.

(v)\ \
In the same vein, it is tempting to mimick example
\ref{ex_cat-cats}(iii), in order to define a full
Godement product $\alpha'*\alpha$ of any two
lax-natural transformations as in (iv) : namely,
one could for example declare this product to be
$(G'*\alpha)\odot(\alpha'*F)$. However, notice that
this lax-natural transformation shall be, in general,
different from $(\alpha'*G)\odot(F'*\alpha)$, which
would be another plausible definition. More precisely,
we have a $2$-cell
$$
\tau^{\alpha'}_{\alpha_A}:
((G'*\alpha)\odot(\alpha'*F))_A\Rightarrow
((\alpha'*G)\odot(F'*\alpha))_A
\qquad
\text{for every $A\in\Ob(\cA)$}
$$
where $\tau^{\alpha'}_{\bullet\bullet}$ denotes the coherence
constraint of $\alpha'$. Thus, the two possible definition
will in general coincide only in case $\alpha'$ is a
strict pseudo-natural transformation.

(vi)\ \
Taking into account remark \ref{rem_equiv-2-cat}(i)
and the observations in (iv) and (v), we conclude
that the datum of all small $2$-categories, all
pseudo-functors between such $2$-categories, and all
lax-natural (or pseudo-natural) transformations, does
{\em not} form a $2$-category, at least not with Godement
products defined as in the foregoing. Rather, the
pseudo-natural isomorphisms exhibited in (iv) hint
at a higher order structure, whose full explication
would require the introduction of a $3$-categorical
formalism that lies beyond the bounds of our treatise.
Instead, we will venture only as far as needed to reach
the notion of {\em modification} of pseudo-natural
transformations (see definition \ref{def_modification}),
which will leave us just outside the threshold of
$3$-category theory, but still will provide a workable
language for expressing the foregoing higher order
compatibilities, and others as well of a similar
nature, that will be encountered in this section.

(vii)\ \
Let us also mention that the datum of all small
$2$-categories, all strict pseudo-functors between
such categories, and all strict pseudo-natural transformations
of such pseudo-functors, does carry a natural $2$-category
structure, with the composition laws for pseudo-functors
and pseudo-natural transformations defined above : the
easy verification shall be left to the reader.
\end{remark}

\begin{example}\label{ex_from-Cats-to-PsFuns}
(i)\ \
Consider two pseudo-functors $F,G:\cA\to\cB$ between
$2$-categories, and $\beta:F\Rightarrow G$ a pseudo-natural
transformation. Denote by $(\gamma^F,\delta^F)$ (resp.
$(\gamma^G,\delta^G)$, resp. $\tau^\beta$) the coherence
constraint of $F$ (resp. $G$, resp. $\beta$).
We deduce from $\beta$ pseudo-natural transformations
$$
\beta^o:G^o\Rightarrow F^o
\qquad
{}^o\beta:{}^oF\Rightarrow{}^oG
$$
with coherence constraints given by the rules
$$
\begin{aligned}
(f^o:A^o\to A'^o)&\,\mapsto\tau^{\beta^o}_{f^o}:=(\tau^{\beta\ -1}_f)^o:
F^of^o\circ\beta^o_{A'^o}\Rightarrow\beta^o_{A^o}\circ G^of^o \\
({}^of:{}^oA\to {}^oA')&\,\mapsto\tau^{{}^o\beta}_{f^o}:=
{}^o(\tau^{\beta\ -1}_f):{}^oG{}^of\circ{}^o\beta_{{}^oA}\Rightarrow
{}^o\beta_{{}^oA'}\circ{}^oF{}^of. 
\end{aligned}
$$
Indeed, the coherence axioms for $\tau^\beta$ imply
the identities
$$
(\tau^{\beta\ -1}_{\one_A})^o\boxminus(\delta^{G\ -1}_A)^o=
(\delta^{F\ -1}_A)^o\boxminus\one_{\beta_A}^o
\qquad
(\gamma^{F\ -1}_{f,f'})^o\boxminus
((\tau^{\beta\ -1}_f)^o\boxvert(\tau^{\beta\ -1}_{f'})^o)=
(\tau^{\beta\ -1}_{f'\circ f})^o\boxminus(\gamma^{G\ -1}_{f,f'})^o
$$
whence :
$$
(\delta^F_A)^o\boxminus(\tau^{\beta\ -1}_{\one_A})^o=
\one_{\beta_A}^o\boxminus(\delta^{G\ -1}_A)^o
\qquad
((\tau^{\beta\ -1}_f)^o\boxvert(\tau^{\beta\ -1}_{f'})^o)
\boxminus(\gamma^G_{f,f'})^o=
(\gamma^F_{f,f'})^o\boxminus(\tau^{\beta\ -1}_{f'\circ f})^o
$$
which in view, of remark \ref{rem_pseudo-funct}(iv),
shows that $\tau^{\beta^o}_\bullet$ fulfills the required
coherence axioms. Likewise one verifies easily the
naturality condition, whence the assertion. The 
corresponding verifications for ${}^o\beta$ are similar,
and shall be left to the reader.

(ii)\ \
By virtue of (i), one can say that -- just as for usual
natural transformation between functors -- the datum of
a pseudo-natural transformation $\beta:F\Rightarrow G$
is equivalent to that of a pseudo-natural transformation
$\beta^o:F^o\Rightarrow G^o$, and also to that of a
pseudo-natural transformation
${}^o\beta:{}^oF\Rightarrow{}^oG$. On the other hand, these
equivalences {\em does not extend to general lax-natural
transformations}. By inspecting the definition,
we see that a lax-natural transformation
$\alpha:{}^oF\Rightarrow{}^oG$ is the datum of a system
of $1$-cells $(\alpha_A:FA\to GA~|~A\in\Ob(\cA))$
together with oriented squares as in remark
\ref{rem_pseudo-natural}(i), but {\em whose orientation
is reversed}, hence $\tau_f$ is a $2$-cell
$\alpha_B\circ Ff\Rightarrow Gf\circ\alpha_A$ in $\cB$.
The naturality and coherence conditions for such a
system of oriented squares must be likewise suitably
reoriented : the details shall be left to the reader.
In some literature, such data
$(\alpha_\bullet,\tau_{\bullet\bullet})$ are called {\em
oplax-natural transformations}, but we shall not use
this terminology.
\end{example}

\begin{example}\label{ex_first-from-Cats-to-PsFuns}
(i)\ \
Let $\cC$ be any $2$-category, and pick a universe $\sU$
such that $\cC$ is $\sU$-small; according to remark
\ref{rem_remark-cat-cats}(i) the product $\cC^o\times\cC$
is representable in $\sU\tdu 2\tdu\bCat$.
Then we get the strict pseudo-functor
$$
H_\cC:\cC^o\times\cC\to\sU\tdu\bCat
$$
(for the natural $2$-category structure on $\sU\tdu\bCat$
discussed in remark \ref{rem_equiv-2-cat}(ii)) that assigns :
\begin{itemize}
\item
to every $C,C'\in\Ob(\cC)$, the category $H_\cC(C,C'):=\cC(C,C')$
\item
to every pair of $1$-cells $f:C_2\to C_1$ and $f':C'_1\to C'_2$
of $\cC$, the functor
$$
H_\cC(f,f'):\cC(C_1,C'_1)\to\cC(C_2,C'_2)
\qquad
g\mapsto f'\circ g\circ f
\quad
\alpha\mapsto f'*\alpha*f 
$$
where $g:C_1\to C'_1$ is any $1$-cell and $\alpha$ any
$2$-cell in $\cC(C_1,C'_1)$
\item
to every pair of $2$-cells $\beta:f\Rightarrow g$,
$\beta':f'\Rightarrow g'$ (with $1$-cells $f,g:C_2\to C_1$
and $f',g':C'_1\to C'_2$), the natural transformation
$$
H_\cC(\beta,\beta'):H_\cC(f,f')\Rightarrow H_\cC(g,g')
\qquad
(h:C_1\to C'_1)\mapsto\beta'*h*\beta.
$$
\end{itemize}

(ii)\ \ 
If $F:\cA\to\cC$ and $G:\cB\to\cC$ are any two pseudo-functors,
we also define
$$
H_\cB(F,G):=H_\cB\circ(F^o\times G):\cA^o\times\cB\to\sU\tdu\bCat.
$$
Likewise, if $F':\cA\to\cC$ and $G':\cB\to\cC$ is another
pair of pseudo-functors, and $\alpha:F\Rightarrow F'$,
$\beta:G\Rightarrow G'$ two pseudo-natural transformations,
we let
$$
H_\cB(\alpha,\beta):=H_\cB*(\alpha^o\times\beta):
H_\cB(F,G)\Rightarrow H_\cB(F',G').
$$
\end{example}

\begin{example}\label{ex_pseudo-functors}
(i)\ \
For every $2$-category $\cA$ and every $X\in\Ob(\cA)$, the
restriction $\ss_X:\cA/X\to\cA$ (resp. $\st_X:X/\cA\to\cA$)
of the pseudo-functor $\ss$ (resp. $\st$) of remark
\ref{rem_remark-cat-cats}(ii) is a well defined strict
pseudo-functor on the $2$-category $\cA/X$ (resp. $X/\cA$
of example \ref{ex_2-arrows}(ii)). Every $1$-cell $h:X\to X'$
induces a strict pseudo-functor
$$
h_*:\cA/X\to\cA/X'
\qquad
(A\xrightarrow{g}X)\mapsto(A\xrightarrow{h\circ g}X').
$$
To every pair of objects $f:A\to X$ and $f':A'\to X$ and
every $1$-cell $(g,\alpha):f\Rightarrow f'$ of $\cA/X$, the
pseudo-functor $h_*$ assigns the $1$-cell
$(g,h*\beta):h_*f\Rightarrow h_*f'$ of $\cA/X'$; lastly,
for every pair of $1$-cells
$(g_1,\alpha_1),(g_2,\alpha_2):f\to f'$ and every $2$-cell
$\beta:g_1\Rightarrow g_2$, we have $h_*\beta:=h*\beta$.

Additionally, every $2$-cell $\beta:h\Rightarrow h'$
between $1$-cells $h,h':X\to X'$ induces a strict
pseudo-natural transformation
$$
\beta_*:h_*\Rightarrow h'_*
$$
which assigns to every $(A\xrightarrow{g}X)\in\Ob(\cA/X)$, the
$1$-cell $(\one_A,\beta*g):h\circ g\to h'\circ g$ of $\cA/X'$.

Lastly, notice that ${}^o({}^o\cA^o/{}^oX^o)^o=X/\cA$; we get
therefore a strict pseudo-functor
$$
h^*:={}^o(({}^oh^o)_*)^o:X'/\cA\to X/\cA
\qquad
(X'\xrightarrow{g}A)\mapsto(X\xrightarrow{g\circ h}A)
$$
and the $2$-cell $\beta$ as in the foregoing induces as
well a strict pseudo-natural transformation
$$
\beta^*:={}^o(({}^o\beta^o)_*)^o:h^*\Rightarrow h'^*.
$$

(ii)\ \
More generally, if $F:\cA\to\cB$ is any pseudo-functor
with coherence constraint $(\delta^F,\gamma^F)$, and
$B\in\Ob(\cB)$ any object, we may define a $2$-category
$$
F\cA/B
$$
whose objects are the pairs $(A,f:FA\to B)$ where $A\in\Ob(\cA)$
and $f$ is a $1$-cell of $\cB$. The $1$-cells
$(g,\alpha):(A,f)\to(A',f')$ consist of a $1$-cell $g:A\to A'$
of $\cA$ and a $2$-cell $\alpha:f\Rightarrow f'\circ Fg$.
For $1$-cells $(g_1,\alpha_1),(g_2,\alpha_2):(A,f)\to(A',f')$,
the $2$-cells $\beta:(g_1,\alpha_1)\Rightarrow(g_2,\alpha_2)$
are the $2$-cells $\beta:g_1\Rightarrow g_2$ such that
$(f'*F\beta)\odot\alpha_1=\alpha_2$. The composition law
for $1$-cells assigns to every pair of $1$-cells
$(A,f)\xrightarrow{(g,\alpha)}(A',f')
\xrightarrow{(g',\alpha')}(A'',f'')$ the $1$-cell
$$
(g',\alpha')\circ(g,\alpha):=
\Bigl(g'\circ g,
\Bigl(\frac{\alpha'}{f'|f''}\boxvert\frac{\alpha}{f|f'}\Bigr)
\boxminus\frac{\gamma^F_{g,g'}}{\one_{FA}|\one_{FA''}}\Bigr).
$$
Let us check the associativity property : let
$(A'',f'')\xrightarrow{(g'',\alpha'')}(A''',f''')$ be
another $1$-cell; then
$$
\begin{aligned}
(g'',\alpha'')\circ((g',\alpha')\circ(g,\alpha))=&\,
(\alpha''\boxvert((\alpha'\boxvert\alpha)\boxminus\gamma^F_{g,g'}))
\boxminus\gamma^F_{g'\circ g,g''} \\
=&\,
((\alpha''\boxminus\one_{Fg''})\boxvert
((\alpha'\boxvert\alpha)\boxminus\gamma^F_{g,g'}))
\boxminus\gamma^F_{g'\circ g,g''} \\
=&\,
(\alpha''\boxvert\alpha'\boxvert\alpha)\boxminus
(\one_{Fg''}\boxvert\gamma^F_{g,g'})\boxminus\gamma^F_{g'\circ g,g''}.
\end{aligned}
$$
On the other hand, a similar calculation yields :
$$
((g'',\alpha'')\circ(g',\alpha'))\circ(g,\alpha)=
(\alpha''\boxvert\alpha'\boxvert\alpha)\boxminus
(\gamma^F_{g',g''}\boxvert\one_{Fg})\boxminus\gamma^F_{g,g''\circ g'}
$$
and the coherence axiom for $F$ translates as :
$(\one_{Fg''}\boxvert\gamma^F_{g,g'})\boxminus
\gamma^F_{g'\circ g,g''}=(\gamma^F_{g',g''}\boxvert\one_{Fg})
\boxminus\gamma^F_{g,g''\circ g'}$, whence the contention.
Then it is easily seen that the identity $1$-cell of any
object $(A,f)$ is given by the pair $(\one_A,\delta^F_A)$ :
details left to the reader. The composition laws for
$2$-cells in $F\cA/B$ are induced by those of $\cA$,
in the obvious way; then the required associativity and
functoriality properties follow straightforwardly. Lastly
we define :
$$
B/F\cA:={}^o({}^o\!F^o({}^o\!\cA^o)/{}^o\!B^o)^o.
$$
The source and target pseudo-functors of (i) generalize as
well : we get strict pseudo-functors
$$
\begin{aligned}
\ss_B&\,:F\cA/B\to\cA
\quad :\quad
(A,f:FA\to B)\mapsto A \\
\st_B&\,:B/F\cA\to\cA
\quad :\quad
(A,f:B\to FA)\mapsto A
\end{aligned}
$$
and every $1$-cell $h:B\to B'$ of $\cB$ induces strict
pseudo-functors
$$
\begin{aligned}
h_*&\,:F\cA/B\to F\cA/B'
\quad :\quad
(A,f:FA\to B)\mapsto(A,h\circ f) \\
h^*&\,:B'/F\cA\to B/F\cA
\quad :\quad
(A,f:B'\to FA)\mapsto(A,f\circ h).
\end{aligned}
$$
Likewise, if $h':B\to B'$ is another $1$-cell of $\cB$, every
$2$-cell $\beta:h'\Rightarrow h'$ of $\cB'$ induces strict
pseudo-natural transformations :
$$
\begin{aligned}
\beta_*&\,:h_*\Rightarrow h'_*
\quad :\quad
(A,f:FA\to B)\mapsto(\one_A,\beta*f) \\
\beta^*&\,:h^*\Rightarrow h'^*
\quad :\quad
(A,f:B\to FA)\mapsto(\one_A,f*\beta).
\end{aligned}
$$ 
\end{example}

\begin{definition}\label{def_modification}
Consider two pseudo-functors $F,G:\cA\to\cB$ between
$2$-categories $\cA$, $\cB$, and two lax-natural
transformations $\alpha,\beta:F\Rightarrow G$.
A {\em modification\/} $\Xi:\alpha\leadsto\beta$ is
a family
$$
\Xi_A:\alpha_A\Rightarrow\beta_A
\qquad
\text{for every $A\in\Ob(\cA)$}
$$
of $2$-cells of $\cB$ that satisfies the
{\em compatibility condition} :
$$
(\Xi_{A'}*Ff)\odot\tau^\alpha_f=\tau^\beta_f\odot(Gf*\Xi_A)
\qquad
\text{every $1$-cell $f:A\to A'$ of $\cA$}
$$
where $\tau^\alpha_{\bullet\bullet}$ (resp.
$\tau^\beta_{\bullet\bullet}$) denotes the coherence
constraint for $\alpha$ (resp. for $\beta$).
\end{definition}

\begin{remark}\label{rem_construct-modifs}
In the situation of definition \ref{def_modification},
let $\Xi:\alpha\leadsto\beta$ be any modification.

(i)\ \
For any pair of $1$-cells $f,g:A\to A'$ of $\cA$, and
any $2$-cell $\gamma:f\Rightarrow g$ we have the identity:
$$
(\Xi_{A'}*F\gamma)\odot\tau^\alpha_f=
\tau^\beta_g\odot(G\gamma*\Xi_A).
$$
Indeed, we may compute :
$$
\begin{aligned}
(\Xi_{A'}*F\gamma)\odot\tau^\alpha_f
=\, &
(\beta_{A'}*F\gamma)\odot(\Xi_{A'}*Ff)\odot\tau^\alpha_f \\
=\, &
(\beta_{A'}*F\gamma)\odot\tau^\beta_f\odot(Gf*\Xi_A) \\
=\, &
\tau^\beta_g\odot(G\gamma*\beta_A)\odot(Gf*\Xi_A) \\
=\, &
\tau^\beta_g\odot(G\gamma*\Xi_A)
\end{aligned}
$$
where the first and the last identities follow from
remark \ref{rem_equiv-2-cat}(i), and the third one
from remark \ref{rem_pseudo-natural}(i). Notice also
that the condition of definition \ref{def_modification}
can be restated as the identity
\set\begin{equation}\label{eq_restate}
\frac{\Xi_{A'}}{\alpha_{A'}|\beta_{A'}}\boxvert
\frac{\tau^\alpha_f}{\alpha_A|\alpha_{A'}}
=\frac{\tau^\beta_f}{\beta_A|\beta_{A'}}\boxvert
\frac{\Xi_A}{\alpha_A|\beta_A}.
\end{equation}

(ii)\ \
Let $\tilde\alpha,\tilde\beta:\cA\to 2\tdu\sMorph(\cB)$
be the pseudo-functors attached to $\alpha$ and $\beta$
as in remark \ref{rem_pseudo-natural}(ii). By unwinding
the definitions, it is easily seen that $\Xi$ is equivalent
to the datum of a strict pseudo-natural transformation
$$
\tilde\Xi:\tilde\alpha\Rightarrow\tilde\beta
\qquad
\text{such that\ \ $\ss*\tilde\Xi=\one_F$\ \ and\ \ 
$\st*\tilde\Xi=\one_G$}.
$$
Namely, we obtain such a $\tilde\Xi$ by assigning, to
every object $A$ of $\cA$, the $1$-cell of $2\tdu\sMorph(\cB)$ :
$$
\xymatrix{
FA \ar[r]^-{\alpha_A} \ddouble 
\drtwocell\omit{_\ \ \ \Xi_A} & GA \ddouble \\
FA \ar[r]_-{\beta_A} & GA.
}$$

(iii)\ \
If $\gamma:F\Rightarrow G$ is any other lax-natural
transformation, and $\Theta:\beta\leadsto\gamma$ any
other modification, we can define the composition 
$$
\Theta\odot\Xi:\alpha\leadsto\gamma
\quad :\quad A\mapsto\Theta_A\odot\Xi_A
\quad\text{for every $A\in\Ob(\cA)$.}
$$
In order to check that $\Theta\odot\Xi$ is indeed a
well defined modification, we may notice that if
$\tilde\Xi$ and $\tilde\Theta$ are defined as in (ii),
then $\Xi\odot\Theta$ is the modification corresponding
to the pseudo-natural transformation
$$
\tilde\Theta\odot\tilde\Xi:\cA\to 2\tdu\sMorph(\cB)
$$
with the composition rule $\odot$ for pseudo-natural
transformations given by remark \ref{rem_pseudo-natural}(iii).
It is clear that this composition law is associative
for any triple of composable modifications.

(iv)\ \
Likewise, if $H:\cA\to\cB$ is another pseudo-functor,
$\alpha',\beta':G\Rightarrow H$ two other lax-natural
transformations, and $\Xi':\alpha'\leadsto\beta'$ any
modification, we can define the modification
$$
\Xi'*\Xi:\alpha'\odot\alpha\leadsto\beta'\odot\beta
\qquad
A\mapsto\Xi'_A*\Xi_A
\qquad\text{for every $A\in\Ob(\cA)$}.
$$
To see that $\Xi'*\Xi$ is well defined, consider -- for
any $1$-cell $f:A\to A'$ in $\cA$ -- the diagram
$$
\xymatrix@+10pt{ FA \ar[r]^-{\alpha_A} \ar[d]_{Ff}
\drtwocell\omit{_\ \ \ \tau^\alpha_f\ } & 
\drtwocell\omit{_\ \ \ \tau^{\alpha'}_f\ }
GA \ar[r]^-{\alpha'_A} \ar[d]|{Gf} & HA \ar[d]^{Hf} \\
FA' \ar[r]|{\alpha_{A'}} \ddouble 
\drtwocell\omit{_\ \ \ \ \ \Xi_{A'}\ } &
\drtwocell\omit{_\ \ \ \ \ \Xi'_{A'}\ } 
GA' \ar[r]|{\alpha'_{A'}} \ddouble & HA' \ddouble \\
FA' \ar[r]_-{\beta_{A'}} & GA' \ar[r]_-{\beta'_{A'}} & HA'.
}$$
We compute :
$$
\begin{aligned}
(\Xi'_{A'}*\Xi_{A'}*Ff)\odot\tau^{\alpha'\circ\alpha}_f
=\, & (\Xi'_{A'}\boxminus\Xi_{A'})\boxvert
(\tau^{\alpha'}_f\boxminus\tau_f^\alpha) \\
=\, & (\Xi'_{A'}\boxvert\tau^{\alpha'}_f)\boxminus
(\Xi_{A'}\boxvert\tau_f^\alpha) \\
=\, & (\tau^{\beta'}_f\boxvert\Xi'_A)\boxminus
(\tau_f^\beta\boxvert\Xi_A) \\
=\, & (\tau^{\beta'}_f\boxminus\tau_f^\beta)\boxvert
(\Xi'_A\boxminus\Xi_A) \\
=\, & \tau_f^{\beta'\circ\beta}\odot(Hf*\Xi'_A*\Xi_A)
\end{aligned}
$$
where the second and fourth equalities follow from
proposition \ref{prop_square-algebra}(ii), and the third
follows from \eqref{eq_restate}. The contention follows.
This operation is clearly associative. As usual, we write
$\alpha*\Xi$ and $\Xi*{\alpha'}$ instead of $\one_\alpha*\Xi$
and respectively $\Xi*\one_{\alpha'}$, for any lax-natural
transformations $\alpha:G\Rightarrow H$ and
$\alpha':E\Rightarrow F$.

(v)\ \
There is a third type of operation for modifications : namely,
let $\cC$ be another $2$-category and $I:\cC\to\cA$ any
pseudo-functor; then we obtain the modification
$$
\Xi\circ I:\alpha*I\Rightarrow\beta*I
\qquad
C\mapsto\Xi_{IC}
\qquad
\text{for every $C\in\Ob(\cC)$}.
$$
To check that $\Xi\circ I$ is indeed a modification, it
suffices to notice that the coherence constraint of
$\alpha*I$ is given by the system of $2$-cells
$\tau^\alpha_{If}$, for $f$ ranging over the $1$-cells
of $\cC$, and correspondingly for the coherence constraint
of $\beta*I$.

(vi)\ \
Likewise, if $K:\cB\to\cC$ is any pseudo-functor, we get
the modification
$$
K\circ\Xi:K*\alpha\leadsto K*\beta
\qquad
A\mapsto K(\Xi_A)
\qquad
\text{for every $A\in\Ob(\cA)$}.
$$
Indeed, notice that the coherence constraints of $K*\alpha$
and $K*\beta$ attach to every $1$-cell $f:A\to A'$ in $\cA$
the $2$-cells
$$
\tau^{K*\alpha}_f:=\gamma^{K\ -1}_{Ff,\alpha_{A'}}\odot K(\tau^\alpha_f)
\odot\gamma^K_{\alpha_A,Gf}
\qquad
\tau^{K*\beta}_f:=\gamma^{K\ -1}_{Ff,\beta_{A'}}\odot K(\tau^\beta_f)
\odot\gamma^K_{\beta_A,Gf}.
$$
Hence we may compute :
$$
\begin{aligned}
(K(\Xi_{A'})*KFf)\odot\tau^{K*\alpha}_f=\,&
(K(\Xi_{A'})*KFf)\odot\gamma^{K\ -1}_{Ff,\alpha_{A'}}\odot
K(\tau^\alpha_f)\odot\gamma^K_{\alpha_A,Gf} \\
=\,&
\gamma^{K\ -1}_{Ff,\beta_{A'}}\odot K(\Xi_{A'}*Ff)\odot
K(\tau^\alpha_f)\odot\gamma^K_{\alpha_A,Gf} \\
=\,&
\gamma^{K\ -1}_{Ff,\beta_{A'}}\odot K((\Xi_{A'}*Ff)\odot\tau^\alpha_f)
\odot\gamma^K_{\alpha_A,Gf} \\
=\,&
\gamma^{K\ -1}_{Ff,\beta_{A'}}\odot K(\tau^\beta_f\odot(\Xi_A*Gf))
\odot\gamma^K_{\alpha_A,Gf} \\
=\,&
\gamma^{K\ -1}_{Ff,\beta_{A'}}\odot K(\tau^\beta_f)\odot K(\Xi_A*Gf)
\odot\gamma^K_{\alpha_A,Gf} \\
=\,&
\gamma^{K\ -1}_{Ff,\beta_{A'}}\odot K(\tau^\beta_f)
\odot\gamma^K_{\beta_A,Gf}\odot(K(\Xi_A)*KGf) \\
=\,&
\tau^{K*\beta}_f\odot(K(\Xi_A)*KGf)
\end{aligned}
$$
where the second and sixth equality follows from remark
\ref{rem_pseudo-funct}(ii), and the third and fifth follow
from remark \ref{rem_pseudo-funct}(i).

(vii)\ \
Let $\Theta$ and $\Xi$ be as in (iii), and $I$ and $K$
like in (v) and (vi); then it is easily seen that
$$
(\Theta\odot\Xi)\circ I=(\Theta\circ I)\odot(\Xi\circ I)
\qquad
K\circ(\Theta\odot\Xi)=(K\circ\Theta)\odot(K\circ\Xi).
$$
Indeed, the first identity is immediate, and the second
follows directly from remark \ref{rem_pseudo-funct}(i).
Furthermore, if $\Xi$ and $\Xi'$ are as in (iv), we have
the identities
$$
(\Xi'*\Xi)\circ I=(\Xi'\circ I)*(\Xi\circ I)
\qquad
K\circ(\Xi'*\Xi)\odot\gamma^K_{\alpha,\alpha'}=
\gamma^K_{\beta,\beta'}\odot((K\circ\Xi')*(K\circ\Xi)).
$$
Indeed, the first identity is immediate, and the second
follows directly from remark \ref{rem_pseudo-funct}(ii).
\end{remark}

\begin{definition}\label{def_PsNat}
Let $\cA$ and $\cB$ be two $2$-categories.

(i)\ \
For every pair of pseudo-functors $F,G:\cA\to\cB$ we
have a category
$$
\sPsNat(F,G)
$$
whose objects are the pseudo-natural transformations
$F\Rightarrow G$, and whose morphisms are the modifications
$\alpha\leadsto\beta$ between them, with the composition rule
$\odot$ given by remark \ref{rem_construct-modifs}(iii).

(ii)\ \
The pseudo-functors $\cA\to\cB$ are the objects of a $2$-category :
$$
\sPsFun(\cA,\cB)
$$
whose $1$-cells $F\to G$ are the pseudo-natural transformations
$F\Rightarrow G$ of, and whose $2$-cells are the modifications.
For fixed $F,G:\cA\to\cB$, the category structure on the set of
$1$-cells $F\to G$ is the one of $\sPsNat(F,G)$. The composition
functor
$$
\sPsNat(F,G)\times\sPsNat(G,H)\to\sPsNat(F,H)
$$
assigns, to any two modifications $\Xi:\alpha\leadsto\beta$
and $\Xi':\alpha'\leadsto\beta'$, the modification $\Xi'*\Xi$.
The associativity and unit axioms for this composition
functor follow by a simple inspection.

(iii)\ \
The strict pseudo-functors $\cA\to\cB$ are the objects of
a $2$-category :
$$
\stPsFun(\cA,\cB)
$$
defined as the sub-$2$-category of $\sPsFun(\cA,\cB)$ whose
$1$-cells are the strict pseudo-natural transformations, and
whose $2$-cells are all the modifications.
\end{definition}

\sset\subsubsection{}\label{subsec_opposite-mods}
In the situation of definition \ref{def_modification},
suppose that $\alpha$ and $\beta$ are pseudo-natural;
in view of example \ref{ex_from-Cats-to-PsFuns}(i), the
modification $\Xi$ induces modifications :
$$
\begin{aligned}
\Xi^o&\,:\alpha^o\leadsto\beta^o
\qquad
A^o\mapsto((\Xi_A^o:\alpha^o_{A^o}\Rightarrow\beta^o_{A^o}) \\
{}^o\Xi&\,:{}^o\beta\leadsto{}^o\alpha
\qquad
{}^oA\mapsto({}^o(\Xi_A):{}^o\beta_{{}^oA}\Rightarrow{}^o\alpha_{{}^oA}).
\end{aligned}
$$
The detailed verification shall be left to the reader. With
the notation of definition \ref{def_PsNat}, it follows that
the rules : $(F:\cA\to\cB)\mapsto(F^o:\cA^o\to\cB^o)$,
$(\beta:F\Rightarrow G)\mapsto(\beta^o:F^o\Rightarrow G^o)$
and $(\Xi:\alpha\leadsto\beta)\mapsto(\Xi^o:\alpha^o\leadsto\beta^o)$
yield a strict isomorphism of $2$-categories :
$$
\sPsFun(\cA,\cB)^o\isom\sPsFun(\cA^o,\cB^o).
$$
Likewise, the rules $F\mapsto{}^oF$, $\beta\mapsto{}^o\beta$
and $\Xi\mapsto{}^o\Xi$ yield a strict isomorphism of
$2$-categories :
$$
{}^o\sPsFun(\cA,\cB)\isom\sPsFun({}^o\cA,{}^o\cB).
$$

\begin{remark}\label{rem_isom-in-psfun}
(i)\ \
The isomorphisms of the category $\sPsFun(\cA,\cB)$ are
precisely the pseudo-natural isomorphisms of pseudo-functors.
Indeed, if $\alpha:F\Rightarrow G$ is such a pseudo-natural
isomorphism, and $\tau$ is its coherence constraint, we obtain
a pseudo-natural isomorphism $\alpha^{-1}:G\Rightarrow F$
by the rule that assigns to every $A\in\Ob(\cA)$ the
$1$-cell $(\alpha_A)^{-1}:GA\to FA$, and to every $1$-cell
$f:A\to A'$ the oriented square
$$
{\diagram GA \ar[r]^-{\alpha_A^{-1}} \ar[d]_{Gf}
\drtwocell\omit{_\ \ \tau'_f} &
FA \ar[d]^{Ff} \\
GA' \ar[r]_-{\alpha_{A'}^{-1}} & FA'
\enddiagram}
\qquad
\text{where $\tau'_f:=\alpha^{-1}_{A'}*\tau_f^{-1}*\alpha^{-1}_A$}.
$$
Indeed, let us check the coherence axioms for $\alpha^{-1}$;
first, we compute, for every $A\in\Ob(\cA)$ :
$$
\begin{aligned}
\tau'_{\one_A}\odot(\delta^F_A*\alpha^{-1}_A)
=(\alpha^{-1}_A*\tau_{\one_A}^{-1}*\alpha^{-1}_A)
\odot(\delta^F_A*\alpha^{-1}_A)
=&\,((\alpha^{-1}_A*\tau_{\one_A}^{-1})\odot\delta^F_A)*\alpha^{-1}_A \\
=&\,\alpha^{-1}_A*(\tau_{\one_A}^{-1}\odot(\alpha_A*\delta^F_A))*\alpha^{-1}_A \\
=&\,\alpha^{-1}_A*(\delta^F_A*\alpha_A)*\alpha^{-1}_A \\
=&\,\alpha^{-1}_A*\delta^F_A
\end{aligned}
$$
which proves the first axiom. Next, for every pair
of $1$-cells $A\xrightarrow{\ f\ }A'\xrightarrow{\ g\ }A''$
of $\cA$ we have :
$$
\begin{aligned}
(\tau'_g\boxvert\tau'_f)\boxminus\gamma^G_{f,g}
=&\,(\alpha^{-1}_{A''}*\gamma^G_{f,g})\odot
(\alpha^{-1}_{A''}*\tau^{-1}_g*\alpha^{-1}_{A'}*Gf)\odot
(Fg*\alpha^{-1}_{A'}*\tau^{-1}_f*\alpha^{-1}_A) \\
=&\,(\alpha^{-1}_{A''}*\gamma^G_{f,g})\odot
(\alpha^{-1}_{A''}*Gg*\tau^{-1}_f*\alpha^{-1}_A)\odot
(\alpha^{-1}_{A''}*\tau^{-1}_g*Ff*\alpha^{-1}_A) \\
=&\,\alpha^{-1}_{A''}*
(\gamma^G_{f,g}\odot(Gg*\tau^{-1}_f*\alpha^{-1}_A)\odot
(\tau^{-1}_g*Ff*\alpha^{-1}_A)) \\
=&\,\alpha^{-1}_{A''}*
((\gamma^G_{f,g}*\alpha_A)\odot(Gg*\tau^{-1}_f)\odot
(\tau^{-1}_g*Ff))*\alpha^{-1}_A
\end{aligned}
$$
and on the other hand
$$
\gamma^F_{f,g}\boxminus\tau'_{g\circ f}=
(\alpha^{-1}_{A''}*\tau^{-1}_{g\circ f}*\alpha^{-1}_A)
\odot(\gamma^F_{f,g}*\alpha^{-1}_A)=
\alpha^{-1}_{A''}*(\tau^{-1}_{g\circ f}\odot(\alpha_{A''}*\gamma^F_{f,g}))
*\alpha^{-1}_A
$$
so we are reduced to showing that
$$
(\gamma^G_{f,g}*\alpha_A)\odot(Gg*\tau^{-1}_f)\odot
(\tau^{-1}_g*Ff)=\tau^{-1}_{g\circ f}\odot(\alpha_{A''}*\gamma^F_{f,g})
$$
which is equivalent to the second coherence axiom for $\alpha$.
Conversely, it is clear that every isomorphism in
$\sPsFun(\cA,\cB)$ must be a pseudo-natural isomorphism.

(ii)\ \
In the same vein, with the notation of definition
\ref{def_PsNat}(i), it is easily seen that the isomorphisms
in the category $\sPsNat(F,G)$ are the modifications
$\Xi:\alpha\leadsto\beta$ such that
$\Xi_A:\alpha_A\Rightarrow\beta_A$ is an isomorphism in
$\cB(FA,GA)$ for every $A\in\Ob(\cA)$ : the details are
left to the reader.
\end{remark}

\begin{example}\label{ex_modifications}
(i)\ \
With the notation of remark \ref{rem_pseudo-natural}(iv),
we may now say that the coherence constraint $\gamma^{G'}$
furnishes an invertible modification
$$
\gamma^{G'}_{\alpha,\beta}:(G'*\beta)\odot(G'*\alpha)\leadsto
G'*(\beta\odot\alpha).
$$

(ii)\ \
Likewise, let us check that the $2$-cell of remark
\ref{rem_pseudo-natural}(v) yields a modification
$$
\tau^{\alpha'}_\alpha:(G'*\alpha)\odot(\alpha'*F)\leadsto
(\alpha'*G)\odot(F'*\alpha)
$$
which is invertible if $\alpha'$ is pseudo-natural.
To this aim, for any $1$-cell $f:A\to A'$ in $\cA$
we consider the diagrams :
$$
\xymatrix@R+10pt@C+18pt
{ F'FA \ar[r]^-{\alpha'_{\!FA}} \ar[d]_{F'Ff}
\drtwocell\omit{_\ \ \ \ \tau^{\alpha'}_{Ff}\ } & 
\drtwocell\omit{_\ \ \ \ \ \ \ \ \ (\tau^\alpha_f)^{G'}\ }
G'FA \ar[r]^-{G'\alpha_{\!A}} \ar[d]|{G'Ff} & G'GA \ar[d]^{G'Gf}
\\
F'FA' \ar[r]|{\alpha'_{\!FA'}} \ar[d]_{F'\alpha_{\!A'}} 
\drtwocell\omit{_\ \ \ \ \ \ \tau^{\alpha'}_{\alpha_{\!A'}}\ } &
\drtwocell\omit{_\ \ \ \ \ \ \ \ \ \ \one_{G'\alpha_{\!A'}}\ } 
G'FA' \ar[r]|{G'\alpha_{A'}} \ar[d]|{G'(\alpha_{\!A'})} &
G'GA' \ddouble
\\
F'GA' \ar[r]_-{\alpha'_{GA'}} & G'GA' \rdouble & G'GA'
}
\quad
\xymatrix@R+10pt@C+18pt
{F'FA \rdouble \ddouble
\drtwocell\omit{_\ \ \ \ \ \ \ \ \one_{F'\alpha_{\!A}}\ } &
F'FA \ar[r]^-{\alpha'_{FA}} \ar[d]|{F'\alpha_A}
\drtwocell\omit{_\ \ \ \ \tau^{\alpha'}_{\alpha_{\!A}}\ } &
G'FA \ar[d]^{G'\alpha_A} \\
F'FA \ar[r]|{F'\alpha_A} \ar[d]_{F'Ff}
\drtwocell\omit{_\ \ \ \ \ \ \ \ (\tau^\alpha_f)^{F'}\ } &
F'GA \ar[r]|{\alpha'_{GA}} \ar[d]|{F'Gf}
\drtwocell\omit{_\ \ \ \ \tau^{\alpha'}_{Gf}\ } &
G'GA \ar[d]^{G'Gf} \\
F'FA' \ar[r]_-{F'\alpha_{A'}} & F'GA' \ar[r]_-{\alpha'_{GA'}}
& G'GA'
}$$
(notation of example \ref{ex_2-funct-of-arrows}(i)) and notice that
$$
\begin{aligned}
(\tau^{\alpha'}_{\alpha_{A'}}*F'Ff)\odot\tau_f^{(G'*\alpha)\odot(\alpha'*F)}
=\lambda:= \, &
(\one_{G'\alpha_{A'}}\boxminus\tau_{\alpha_{A'}}^{\alpha'})\boxvert
((\tau^\alpha_f)^{G'}\boxminus\tau^{\alpha'}_{Ff}) \\
=\, &
(\one_{G'\alpha_{A'}}\boxvert(\tau^\alpha_f)^{G'})\boxminus
(\tau_{\alpha_{A'}}^{\alpha'}\boxvert\tau^{\alpha'}_{Ff}) \\
\tau^{(\alpha'*G)\odot(F'*\alpha)}_f\odot(G'Gf*\tau^{\alpha'}_{\alpha_A})
=\mu:= \, & 
(\tau^{\alpha'}_{Gf}\boxminus(\tau^{\alpha'}_f)^{F'})\boxvert
(\tau^{\alpha'}_{\alpha_A}\boxminus\one_{F'\alpha_A}) \\
=\, &
(\tau^{\alpha'}_{Gf}\boxvert\tau^{\alpha'}_{\alpha_A})\boxminus
((\tau^{\alpha'}_f)^{F'}\boxvert\one_{F'\alpha_A}).
\end{aligned}
$$
So we have only to check that
$$
\lambda\boxminus\frac{\one_{F'\alpha_{A'}\circ F'Ff}}{\one_{F'FA}|F'\alpha_{A'}}
=\frac{\one_{G'Gf\circ G'\alpha_A}}{G'\alpha_A|\one_{G'GA}}\boxminus\mu.
$$
However, the definition of $2\tdu\sMorph(G')$ and the
coherence axiom for $\alpha'$ yield respectively
$$
\begin{aligned}
\one_{G'\alpha_{A'}}\boxvert(\tau^\alpha_f)^{G'}
=\, &
\frac{\gamma^{G'}_{\alpha_A,Gf}}{G'(\alpha_A)|\one_{G'GA'}}
\boxminus
\frac{G'(\tau^\alpha_f)}{\one_{G'FA}|\one_{G'GA'}}
\boxminus
\frac{(\gamma^{G'}_{Ff,\alpha_{A'}})^{-1}}{\one_{G'FA}|\one_{G'GA'}}
 \\
\tau_{\alpha_{A'}}^{\alpha'}\boxvert\tau^{\alpha'}_{Ff}
=\, &
\frac{\gamma^{G'}_{Ff,\alpha_{A'}}}{\one_{G'FA}|\one_{G'GA'}}
\boxminus
\frac{\tau^{\alpha'}_{\alpha_{A'}\circ Ff}}{\alpha'_{FA}|\alpha'_{GA'}}
\boxminus
\frac{(\gamma^{F'}_{Ff,\alpha_{A'}})^{-1}}{\one_{F'FA}|\one_{F'GA'}}
\end{aligned}
$$
therefore :
$$
\begin{aligned}
\lambda=
\,& \gamma^{G'}_{\alpha_A,Gf}\boxminus G'(\tau^\alpha_f)
\boxminus\tau^{\alpha'}_{\alpha_{A'}\circ Ff}\boxminus
(\gamma^{F'}_{Ff,\alpha_{A'}})^{-1} \\
=\,& \gamma^{G'}_{\alpha_A,Gf}\boxminus
\frac{\tau^{\alpha'}_{Gf\circ\alpha_A}}{\alpha'_{FA}|\alpha'_{GA'}}
\boxminus\frac{F'(\tau^\alpha_f)}{\one_{F'FA}|\one_{F'GA'}}
\boxminus(\gamma^{F'}_{Ff,\alpha_{A'}})^{-1} \\
=\,& \frac{\one_{G'Gf\circ G'\alpha_A}}{G'\alpha_A|\one_{G'GA'}}
\boxminus(\tau^{\alpha'}_{Gf}\boxvert\tau^{\alpha'}_{\alpha_A})
\boxminus\frac{(\gamma^{F'}_{\alpha_A,Gf})}{\one_{F'FA}|\one_{F'GA'}}
\boxminus\frac{(\gamma^{F'}_{\alpha_A,Gf})^{-1}}{\one_{F'FA}|\one_{F'GA'}}
\boxminus\frac{(\tau^\alpha_f)^{F'}}{\one_{F'FA}|\one_{F'GA'}} \\
=\,& \frac{\one_{G'Gf\circ G'\alpha_A}}{G'\alpha_A|\one_{G'GA'}}
\boxminus(\tau^{\alpha'}_{Gf}\boxvert\tau^{\alpha'}_{\alpha_A})
\boxminus\frac{(\tau^\alpha_f)^{F'}}{\one_{F'FA}|\one_{F'GA'}}
\end{aligned}
$$
where the second equality follows from remark
\ref{rem_pseudo-natural}(i) and the coherence axiom for
$\alpha'$. The sought identity follows straightforwardly.

(iii)\ \
For every $2$-category $\cA$, we denote by $i_\cA$ the
identity (strict) pseudo-natural transformation of $\one_\cA$.
In other words, $i_\cA$ assigns to every $A\in\Ob(\cA)$ the
$1$-cell $\one_A$. Let $\cB,\cC$ be any other $2$-categories,
and $F:\cA\to\cB$, $G:\cC\to\cA$ any two pseudo-functors;
clearly
$$
i_\cA*G=\one_G.
$$
On the other hand, $F*i_\cA$ is a pseudo-natural transformation
that assigns to every $A\in\Ob(\cA)$ the $1$-cell $F\one_A:FA\to FA$.
If $(\delta^F,\gamma^F)$ is the coherence constraint for $F$,
the coherence constraint of $F*i_\cA$ assigns to every $1$-cell
$f:A\to A'$ of $\cA$ the $2$-cell
$\one^F_f:=\gamma^{F\ -1}_{f,\one_{A'}}\odot\gamma^F_{\one_A,f}$. Then,
the unit axiom for $F$ implies that the system of $2$-cells
$(\delta^F_A~|~A\in\Ob(\cA))$ defines an invertible modification
$$
\delta^F:\one_F\leadsto F*i_\cA.
$$
\end{example}

\begin{remark}\label{rem_reduced-2-cats}
As already pointed out (remark \ref{rem_pseudo-natural}(vi)),
the datum of the category $2\tdu\bCat$ and the system of
categories $(\sPsFun(\cA,\cB)~|~\cA,\cB\in\Ob(2\tdu\bCat))$
does not amount to a $2$-category. However, example
\ref{ex_modifications} suggests a straightforward
variation that does result in an interesting $2$-category.
Namely, for every pair of categories $\cA$ and $\cB$, let
us denote by
$$
\overline{\sPsFun}(\cA,\cB)
$$
the category whose objects are the pseudo-functors
$\cA\to\cB$ and whose morphisms are the equivalence
classes $[\beta]$ of pseudo-natural transformations
$\beta$ between such pseudo-functors, where we declare
that two pseudo-natural transfomations
$\alpha,\beta:F\Rightarrow G$ are equivalent if and only
if there exists an invertible modification
$\alpha\leadsto\beta$. The composition law for such
equivalence classes is given by the rule
$$
[\beta]\odot[\beta']:=[\beta\odot\beta']
$$
for every two pseudo-natural transformations
$\beta:F\Rightarrow G$, $\beta':G\Rightarrow H$ of
pseudo-functors $F,G,H:\cA\to\cB$. Then, if $\cC$
is another $2$-category, $F',G':\cB\to\cC$ and
$\alpha:F'\Rightarrow G'$ another pseudo-natural
transformation, we set
$$
[\alpha]*[\beta]:=[(G'*\beta)\odot(\alpha*F)].
$$
By virtue of example \ref{ex_modifications}, this operation
yields then a well defined functor
$$
\overline\sPsFun(\cA,\cB)\times\overline\sPsFun(\cB,\cC)\to
\overline\sPsFun(\cA,\cC)
$$
that satisfies the associativity axiom for the composition
bifunctor in a $2$-category (see \eqref{sec_2Cats}).
Thus, for every universe $\sU'$, the $\sU'$-small $2$-categories
are the objects of a $2$-category
$$
\sU'\tdu\overline{2\tdu\bCat}
$$
that we call the {\em reduced $2$-category of $\sU'$-small
$2$-categories}, whose $1$-cells are the pseudo-functors and
whose $\Hom$-categories are the foregoing categories
$\overline\sPsFun(-,-)$. As usual, when $\sU'=\sU$ we write
simply $\overline{2\tdu\bCat}$ for this $2$-category.
\end{remark}

\begin{remark}\label{rem_remains-of-the-day}
(i)\ \
Let $\cA$ and $\cB$ be two $2$-categories. Every pair of
pseudo-functors $F:\cA'\to\cA$ and $G:\cB\to\cB'$ induces
a pseudo-functor
$$
\sPsFun(F,G):\sPsFun(\cA,\cB)\to\sPsFun(\cA',\cB')
\qquad
(H:\cA\to\cB)\mapsto G\circ H\circ F
$$
that assigns to every pseudo-natural transformation
$\beta:H\Rightarrow H'$ between pseudo-functors
$H,H':\cA\to\cB$ the pseudo-natural transformation
$G*\beta*F:G\circ H\circ F\Rightarrow G\circ H'\circ F$,
and to every modification $\Xi:\beta\leadsto\beta'$
between such pseudo-natural transformations $\beta,\beta'$,
the modification $G\circ\Xi\circ F:G*\beta*F\leadsto G*\beta'*F$.
The coherence constraint of $\sPsFun(F,G)$ is the pair
$(\delta^{(F,G)}_\bullet,\gamma^{(F,G)}_{\bullet\bullet})$ defined
as follows :
$$
\delta^{(F,G)}_H:=\delta^G\circ H\circ F:\one_{GHF}\leadsto G*\one_{HF}
\qquad
\text{for every pseudo-functor $H:\cA\to\cB$}
$$
where $\delta^G:\one_G\leadsto G*i_\cA$ is the invertible
modification provided by example \ref{ex_modifications}(iii), and
$$
\gamma^{(F,G)}_{\beta,\beta'}:=\gamma^G_{\beta,\beta'}\circ F:
(G*\beta'*F)\odot(G*\beta*F)\leadsto G*(\beta'\odot\beta)*F
$$
for every pair of pseudo-natural transformations
$\beta:H\Rightarrow H'$ and $\beta':H'\Rightarrow H'$,
where $\gamma^G{\beta,\beta'}$ is the invertible modification
as in example \ref{ex_modifications}(i).
The verification is straightforward, in view of remark
\ref{rem_construct-modifs}(vii). In case $F=\one_\cA$
(resp. $G=\cB$), we also write $\sPsFun(\cA,G)$ (resp.
$\sPsFun(F,\cB)$) for this pseudo-functor. Notice that
$\sPsFun(F,\cB)$ is strict for every pseudo-functor $F$,
whereas $\sPsFun(\cA,G)$ is strict only if $G$ is strict.

(ii)\ \
Every pseudo-natural transformation $\alpha:F\Rightarrow F'$
between pseudo-functors $F,F':\cA'\to\cA$ induces a
pseudo-natural transformation
$$
\sPsFun(\alpha,\cB):\sPsFun(F,\cB)\Rightarrow\sPsFun(F',\cB)
\qquad
(H:\cA\to\cB)\mapsto(H*\alpha:H\circ F\Rightarrow H\circ F')
$$
whose coherence constraint assigns to every pseudo-natural
transformation $\beta:H\Rightarrow H'$ the oriented square :
$$
\xymatrix@C+20pt{
H\circ F \ar@{=>}[r]^-{H*\alpha} \ar@{=>}[d]_{\beta*F}
\drtwocell\omit{_\ \ \ \ \ \ \ (\tau^\beta_\alpha)^{-1}} &
H\circ F' \ar@{=>}[d]^{\beta*F'} \\
H'\circ F \ar@{=>}[r]_-{H'*\alpha} & H'\circ F'
}$$
where $\tau^\beta_\alpha$ is the invertible modification
provided by example \ref{ex_modifications}(ii). Indeed,
the naturality of the rule $\beta\mapsto(\tau^\beta_\alpha)^{-1}$
with respect to modifications $\Xi:\beta\leadsto\beta'$
follows directly from the compatibility condition for $\Xi$.
To check the coherence axioms, notice first that
$\tau^{\one_H}_\alpha=\one_{H*\alpha}$ for every pseudo-functor
$H:\cA\Rightarrow\cB$. Lastly, if $\beta:H\Rightarrow H'$
and $\beta':H'\Rightarrow H''$ are two pseudo-natural
transformations, we need to show the identity
$$
(\tau^{\beta'}_\alpha)^{-1}\boxvert(\tau^\beta_\alpha)^{-1}=
(\tau^{\beta'\odot\beta}_\alpha)^{-1}
$$
or equivalently :
$\tau^{\beta'}_\alpha\boxminus\tau^\beta_\alpha=\tau^{\beta'\odot\beta}_\alpha$,
which follows by direct inspection.

(iii)\ \
Likewise, every pseudo-natural transformation
$\lambda:G\Rightarrow G'$ between pseudo-functors
$G,G':\cB\to\cB'$ induces a pseudo-natural transformation
$$
\sPsFun(\cA,\lambda):\sPsFun(\cA,G)\Rightarrow\sPsFun(\cA,G')
\qquad
(H:\cA\to\cB)\mapsto(\lambda*H:G\circ F\Rightarrow G'\circ F)
$$
whose coherence constraint assigns to every pseudo-natural
transformation $\beta:H\Rightarrow H'$ the oriented square
$\tau^\lambda_\beta$ provided by example \ref{ex_modifications}(ii).
Indeed, the naturality of the rule : $\beta\mapsto\tau^\lambda_\beta$
with respect to modifications $\Xi:\beta\leadsto\beta'$ follows
from the naturality of $\tau^\lambda$, applied to the $2$-cells
$\Xi_A:\beta_A\to\beta'_A$, for every $A\in\Ob(\cA)$. The
coherence condition for $\tau^\lambda_{\one_H}$ follows from
the coherence condition for $\tau^\lambda_{HA}$, for every
pseudo-functor $H:\cA\to\cB$ and every $A\in\Ob(\cA)$.
Likewise, the coherence condition relative to a composable
pair of pseudo-natural transformations $\beta:H\Rightarrow H'$,
$\beta':H'\Rightarrow H''$ follows from the coherence condition
for $\tau^\lambda$, relative to the composable pair of $1$-cells
$\beta_A,\beta'_A$, for every $A\in\Ob(\cA)$.

(iv)\ \
Furthermore, for every $F,F'$ as in (ii), every modification
$\Xi:\alpha\leadsto\beta$ between pseudo-natural transformations
$\alpha,\beta:F\Rightarrow F'$ induces a modification
$$
\sPsFun(\Xi,\cB):\sPsFun(\alpha,\cB)\leadsto\sPsFun(\beta,\cB)
\qquad
(G:\cA\to\cB)\mapsto(G\circ\Xi:G*\alpha\leadsto G*\beta).
$$
Indeed, the compatibility condition for this modification
amounts to the identity :
$$
(G'(\Xi_X)*\gamma_{FX})\odot(\tau^\gamma_{\alpha_X})^{-1}=
(\tau^\gamma_{\beta_X})^{-1}\odot(\gamma_{F'X}*G(\Xi_X))
\qquad
\text{for every $X\in\Ob(\cA)$}
$$
which follows easily from the naturality of $\tau^\gamma$.

(v)\ \
Suppose moreover that $\cB$ is small; then we deduce
a strict pseudo-functor
$$
\sPsFun(-,\cB):\overline{2\tdu\bCat}{}^o\to\overline{2\tdu\bCat}
\qquad
\cA\mapsto\sPsFun(\cA,\cB)
$$
acting on the $1$-cells and $2$-cells of $\overline{2\tdu\bCat}{}^o$
as explicited in (i) and (ii). Indeed, if $F:\cA'\to\cA$ and
$F':\cA''\to\cA'$ are two pseudo-functors between small
$2$-categories, a simple inspection shows that
$\sPsFun(F',\cB)\circ\sPsFun(F,\cB)=\sPsFun(F\circ F',\cB)$.
Moreover, (iv) implies that if $\alpha,\beta:F\Rightarrow F'$
are pseudo-natural transformations with $[\alpha]=[\beta]$
in $\overline{2\tdu\bCat}$, then
$[\sPsFun(\alpha,\cB)]=[\sPsFun(\beta,\cB)]$.
Next, for every three pseudo-functors $F,F',F'':\cA'\to\cA$ and
every pseudo-natural transformations $\alpha:F\Rightarrow F'$
and $\alpha':F'\Rightarrow F''$, with the notation of remark
\ref{rem_reduced-2-cats}, we need to check that :
\set\begin{equation}\label{eq_uffa-che-palle}
[\sPsFun(\alpha',\cB)\odot\sPsFun(\alpha,\cB)]=
[\sPsFun(\alpha'\odot\alpha,\cB)].
\end{equation}
Now, to every pseudo-functor $G:\cA\to\cB$, the left-hand side
of \eqref{eq_uffa-che-palle} assigns the class of
$(G*\alpha')\odot(G*\alpha)$, whereas the right-hand side assigns
the class of $G*(\alpha'\odot\alpha)$. These two pseudo-natural
transformations are related by the invertible modification
$\gamma^G_{\alpha,\alpha'}$ of example \ref{ex_modifications}(i),
so we come down to checking that the rule :
$G\mapsto\gamma^G_{\alpha,\alpha'}$ defines an invertible
modification $\sPsFun(\alpha',\cB)\odot\sPsFun(\alpha,\cB)
\leadsto\sPsFun(\alpha'\odot\alpha,\cB)$. The latter assertion
is an immediate consequence of the coherence axioms for
$\tau^\beta$ : details left to the reader. Lastly, let us
check that for every two pairs of pseudo-functors
$F_1,F'_1:\cA'\to\cA$, $F_2,F'_2:\cA''\to\cA'$ and every
pseudo-natural transformations $\alpha_1:F_1\Rightarrow F'_1$,
$\alpha_2:F_2\Rightarrow F'_2$ we have :
\set\begin{equation}\label{eq_uffa-che-palle-last}
[\sPsFun(\alpha_2,\cB)]*[\sPsFun(\alpha_1,\cB)]=
[\sPsFun([\alpha_1]*[\alpha_2],\cB)].
\end{equation}
We may consider separately the case where $F_2=F'_2$ and
$\alpha_2=\one_{F_2}$, and the case where $F_1=F'_1$ and
$\alpha_1=\one_{F_1}$. In the first case, the left-hand side
assigns to every $G:\cA\to\cB$ the pseudo-natural
transformation $(G*\alpha_1)*F_2$, whereas the right-hand
side assigns $G*(\alpha_1)*F_2$; but these coincide, and
likewise, the coherence constraints of both sides is given
by the rule : $\beta\mapsto(\tau^\beta_{\alpha_1*F_2})^{-1}$.
One argues similarly in the second case : details left to
the reader.
\end{remark}

\subsection{The formalism of base change}
\label{sec_link}
Let $\cA$ be a $2$-category, and $A,B$ two objects of $\cA$.
A {\em link from $A$ to $B$ in $\cA$} is a datum
$$
\cL:=(F,G,\eta,\eps):A\to B
$$
consisting of :
\begin{itemize}
\item
$1$-cells $F:A\to B$ and $G:B\to A$
\item
$2$-cells $\eta:\one_B\Rightarrow F\circ G$
and $\eps:G\circ F\Rightarrow\one_A$
\end{itemize}
such that $\eta$ and $\eps$ are related by the {\em triangular
identities} :
$$
(F*\eps)\odot(\eta*F)=\one_F
\qquad\text{and}\qquad
(\eps*G)\odot(G*\eta)=\one_G.
$$
Especially, $(G,F)$ is an adjoint pair of $1$-cells of $\cA$.
We say that $\eta$ and $\eps$ are the {\em unit} and {\em counit}
of the link $\cL$. We define a composition law for links as
follows. With $\cL$ as in the foregoing, a third object $C$
of $\cA$ and a second link $\cL':=(F',G',\eta',\eps'):B\to C$
we let
$$
\cL'\circ\cL:=(F'\circ F,G\circ G',\eta^{\cL'\circ\cL},\eps^{\cL'\circ\cL}):
A\to C
$$
where $\eta^{\cL'\circ\cL}:=(F'*\eta*G')\odot\eta'$ and
$\eps^{\cL'\circ\cL}:=\eps\odot(G*\eps'*F)$. We need to check that
$\alpha:=(F'F*(\eps\odot G*\eps'*F))\odot
((F'*\eta*G'\odot\eta')*F'F)=\one_{F'F}$. We compute :
$$
\begin{aligned}
\alpha&\,=
(F'*((F*\eps)\odot(FG*\eps'*F)))\odot(((F'*\eta*G'F')\odot(\eta'*F'))*F) \\
&\,=
(F'*((F*\eps)\odot(\eta*F)))\odot(((F'*\eps')\odot(\eta'*F'))*F) \\
&\,=(F'*\one_F)\odot(\one_{F'}*F) \\
&\,=\one_{F'F}
\end{aligned}
$$
where the second equality holds by remark \ref{rem_equiv-2-cat}(i)
and the third follows from the triangular identities for $(\eta,\eps)$
and $(\eta',\eps')$. Likewise one checks the second triangular
identity for $(\eta'',\eps'')$ (the details are left to the reader).

\begin{remark}\label{rem_compose-links}
(i)\ \
Notice that the composition law for links is associative;
namely, in the situation of \eqref{sec_link}, if
$\cL'':=(F'',G'',\eta'',\eps''):C\to D$ is another link,
we have $(\cL''\circ\cL')\circ\cL=\cL''\circ(\cL'\circ\cL)$.
Indeed, the assertion boils down to the identities :
$$\begin{aligned}
(F''*((F'*\eta*G')\odot\eta')*G'')\odot\eta''&\,=
(F''F'*\eta*G'G'')\odot(F''*\eta'*G'')\odot\eta'' \\
\eps\odot(G*\eps'*F)\odot(GG'*\eps''*F'F)&\,=
\eps\odot(G*(\eps'\odot(G'*\eps''*F'))*F)
\end{aligned}
$$
which are both clear.

(ii)\ \
Notice also that every link $\cL:=(F,G,\eta,\eps):A\to B$
induces two {\em opposite links}
$$
\begin{aligned}
\cL^o&\,:=(G^o,F^o,\eta^o,\eps^o):A^o\to B^o
&\qquad\qquad & \text{in $\cA^o$} \\
{}^o\!\cL&\,:=({}^oG,{}^o\!F,{}^o\eps,{}^o\eta):{}^o\!B\to{}^o\!\!A
&\qquad\qquad & \text{in ${}^o\!\!\cA$}.
\end{aligned}  
$$
\end{remark}

\sset\subsubsection{}\label{subsec_transform-links}
Let $\cL,\cL':A\to B$ be two links, with
$\cL=(F,G,\eta,\eps)$ and $\cL'=(F',G',\eta',\eps')$;
a {\em transformation from $\cL$ to $\cL'$}, denoted
$$
\beta:\cL\Rightarrow\cL'
$$
is defined as a $2$-cell $\beta:F\Rightarrow F'$. The standard
operations on $2$-cells can be upgraded easily to transformations
of links :

$\bullet$\ \
If $\cL'':A\to B$ is another link and $\beta':\cL'\Rightarrow\cL''$
a second transformation, we get the transformation
$$
\beta'\odot\beta:\cL\Rightarrow\cL''.
$$

$\bullet$\ \
For a diagram of links and transformations of links :
\set\begin{equation}\label{eq_Godement-and-links}
{\diagram A \rtwocell^\cL_{\cL'}{\beta} & B
\rtwocell^\cP_{\cP'}{\alpha} & C
\enddiagram}
\end{equation}
we get the transformation
$$
\alpha*\beta:\cP\circ\cL\Rightarrow\cP'\circ\cL'.
$$
For $\alpha=\one_\cP$ (resp. for $\beta=\one_\cL$) we also
denote this transformation by $\cP*\beta$ (resp. by $\alpha*\cL$).
Thus, the links of $\cA$ are the $1$-cells of a $2$-category
$$
\sLink(\cA)
$$
whose objects are the objects of $\cA$ and whose $2$-cells
are the transformations of links, with the compositions laws
for $1$-cells and $2$-cells given in the foregoing.

$\bullet$\ \
For a transformation $\beta:\cL\Rightarrow\cL'$, we define an
{\em adjoint} $2$-cell of $\cA$ :
$$
\beta^\dagger:=(\eps'*G)\odot(G'*\beta*G)\odot(G'*\eta):G'\Rightarrow G.
$$
The operations of remark \ref{rem_compose-links}(ii) can be
combined with this adjoint construction to produce natural
isomorphisms of $2$-categories; indeed, we have :

\begin{proposition}\label{prop_isoms-of-links}
For every $2$-category $\cA$ we have strict isomorphisms
of $2$-categories :
$$
{}^o\sLink(\cA)\isom\sLink(\cA^o)
\qquad
\sLink(\cA)^o\isom\sLink({}^o\!\!\cA).
$$
\end{proposition}
\begin{proof} Indeed, the first strict pseudo-functor is given
by the rules :
$$
A\mapsto A^o
\qquad
(\cL:A\to B)\mapsto(\cL^o:A^o\to B^o)
\qquad
(\beta:\cL\Rightarrow\cL')\mapsto(\beta^{\dagger o}:\cL'^o\Rightarrow\cL^o)
$$
and the second one is given by the rules :
$$
A\mapsto{}^o\!\!A
\qquad
(\cL:A\to B)\mapsto({}^o\!\cL:{}^o\!B\to{}^o\!\!A)
\qquad
(\beta:\cL\Rightarrow\cL')\mapsto
({}^o\!\beta^\dagger:{}^o\!\cL\Rightarrow{}^o\!\cL')
$$
for every object $A$, every $1$-cell $\cL$ and every $2$-cell
$\beta$ of $\sLink(\cA)$. The verification that these rules
yield strict pseudo-functors comes down to the following :

\begin{claim}
(i)\ \
Let $\cL,\cL',\cL'':A\to B$ be three links of $\cA$,
and $\beta:\cL\Rightarrow\cL'$, $\beta':\cL'\Rightarrow\cL''$
two transformations. Then we have
$\beta^\dagger\odot\beta'^\dagger=(\beta'\odot\beta)^\dagger$.

(ii)\ \
In the situation of diagram \eqref{eq_Godement-and-links}
we have $(\alpha*\beta)^\dagger=\beta^\dagger*\alpha^\dagger$.
\end{claim}
\begin{pfclaim}Say that $\cL=(F,G,\eta,\eps)$,
$\cL'=(F',G',\eta',\eps')$, $\cL''=(F'',G'',\eta'',\eps'')$.
Recall that every $1$-cell $f:X\to Y$ of $\cA$ induces functors
$$
f^Z_*:\cA(Z,X)\to\cA(Z,Y)
\qquad
(h:Z\to X)\mapsto f\circ h
\qquad
(\nu:h\Rightarrow h')\mapsto f*\nu
$$
for every $Z\in\Ob(\cA)$, and every $2$-cell
$\lambda:f\Rightarrow g$ of $\cA$ induces natural transformations
$$
\lambda^Z_*:f^Z_*\Rightarrow g^Z_*
\qquad
(h:Z\to X)\mapsto(\lambda*h:f^Z_*(h)\Rightarrow g^Z_*(h))
$$
(see the proof of lemma \ref{lem_long-time}(i)). Then it is
clear that for every $Z\in\Ob(\cA)$, every link
$\cP:=(H,K,\eta^\cP,\eps^\cP)$ induces an adjoint pair of
functors $(H^Z_*,K^Z_*)$ with unit $(\eta^\cP)^Z$ and counit
$(\eps^\cP)^Z$. With this notation, it is easily seen that
\set\begin{equation}\label{eq_setting-back-pg-no}
(\beta^\dagger)^Z_*=(\beta^Z_*,\eta^Z_*,\eta'^Z_*)^\dagger
\end{equation}
where the right-hand side is the adjoint of the natural
transformation $\beta^Z_*$, defined as in remark
\ref{rem_adjoint-transf}(ii) : the details are left to
the reader. Then both (i) and (ii) follow straightforwardly
from remark \ref{rem_adjoint-transf}(iv,v).
\end{pfclaim}

Lastly, in order to see that these pseudo-functors are
isomorphisms of $2$-categories, it suffices to show more
precisely that for every transformation of links
$\beta:\cL\Rightarrow\cL'$ we have :
$$
(\beta^\dagger)^\dagger=\beta.
$$
Again, by virtue of \eqref{eq_setting-back-pg-no}, the
assertion is reduced to the corresponding identity for
adjoints of natural transformations, and the latter is
known by virtue of remark \ref{rem_adjoint-transf}(ii).
\end{proof}

\sset\subsubsection{}\label{subsec_adjunction-from-link}
Notice that every $1$-cell $F:A\to B$ of the $2$-category
$\cA$ induces a strict pseudo-functor
$$
F_*:\cA/A\to\cA/B
\qquad
(g:X\to A)\mapsto(F\circ g:X\to B)
$$
that maps every $1$-cell
$(h,\alpha):(X\xrightarrow{g}A)\to(Y\xrightarrow{g'}A)$ of
$\cA/A$ to the $1$-cell of $\cA/B$
$$
F*(h,\alpha):=(h,F*\alpha)
$$
and every $2$-cell $\beta:(h_1,\alpha_1)\Rightarrow(h_2,\alpha_2)$
of $\cA/A$ to the $2$-cell
$\beta:(h_1,F*\alpha_1)\Rightarrow(h_2,F*\alpha_2)$ of $\cA/B$
(notation of example \ref{ex_2-arrows}(ii)). With this notation,
we claim that every link $\cL:=(F,G,\eta,\eps):A\to B$ of $\cA$
induces an adjunction for the pair of functors
$F_*:\cA/A\to\cA/B$ and $G_*:\cA/B\to\cA/A$. Namely, we
have a natural bijection :
$$
\theta^\cL_{g,f}:
\cA/B(Y\xrightarrow{g}B,X\xrightarrow{F\circ f}B)\isom
\cA/A(Y\xrightarrow{G\circ g}A,X\xrightarrow{f}A)
\quad
(h,\alpha)\mapsto(h,(\eps*f)\odot(G*\alpha))
$$
whose inverse
$\mu^\cL_{g,f}:\cA/A(G\circ g,f)\isom\cA/B(g,F\circ f)$
is given by the rule :
$$
(k,\beta)\mapsto(k,(F*\beta)\odot(\eta*g)).
$$
The verification of the naturality of $\theta^\cL_{f,g}$ with
respect to $f$ and $g$ shall be left to the reader. Let
us check that $\mu^\cL_{g,f}\circ\theta^\cL_{g,f}$ is the
identity map on the set $\cA/B(g,F\circ f)$ : indeed,
for every $1$-cell $(h,\alpha)$, this composition is
computed as the composition of the following oriented
squares :
$$
\xymatrix{ Y \ar[r]^-h \ar[d]_g \ddrtwocell\omit{^\alpha\ } &
X \rdouble \drtwocell\omit{^\one\ } \ar[d]_f & X \ar[d]^f \\
B \ddouble & A \ar[d]_F \rdouble \drtwocell\omit{^\eps\ } &
A \ddouble \\
B \rdouble \ddouble \drtwocell\omit{^\one\ } &
B \ddouble \ar[r]^-G \drtwocell\omit{^\eta\ } & A \ar[d]^F \\
B \rdouble & B \rdouble & B.
}$$
Then the assertion follows easily, by applying the triangular
identities for the pair $(\eta,\eps)$ and invoking as usual
proposition \ref{prop_square-algebra}. Likewise we check
that $\theta^\cL_{g,f}\circ\mu^\cL_{g,f}$ is the identity on
$\cA/A(G\circ g,f)$ : the details shall be left to the reader.

\sset\subsubsection{}\label{subsec_base-change-map}
Consider now an {\em oriented square} $\cD$ of links,
together with its {\em adjoint squares} :
$$
\cD\ :\
{\spreaddiagramcolumns{20pt}
\diagram A' \ar[r]^-{\cL'} \ar[d]_{\cM'}
\drtwocell\omit{_\ \beta} & B' \ar[d]^\cM \\
A \ar[r]_-\cL & B
\enddiagram}
\qquad
\cD^\dagger\ :\
{\spreaddiagramcolumns{20pt}
\diagram A'^o \ar[r]^-{\cL'^o} \ar[d]_{\cM'^o}
\drtwocell\omit{^\ \beta^{\dagger o}\ \ } & B'^o \ar[d]^{\cM^o} \\
A^o \ar[r]_-{\cL^o} & B^o
\enddiagram}
\qquad
{}^\dagger\cD\ :\
{\spreaddiagramcolumns{20pt}
\diagram {}^o\!B \ar[r]^-{{}^o\!\!\cL} \ar[d]_{{}^o\!\!\cM}
\drtwocell\omit{^{}^o\!\beta^\dagger\ \ } &
{}^o\!\!A \ar[d]^{{}^o\!\!\cM'} \\
{}^o\!B' \ar[r]_-{{}^o\cL'} & {}^o\!\!A'
\enddiagram}$$
{\em i.e.} the orientation $\beta:\cM\circ\cL'\Rightarrow\cL\circ\cM'$
is a transformation of links, so that ${}^o\!\beta^\dagger$ is the
adjoint transformation
${}^o\!\!\cL'\circ{}^o\!\!\!\cM\Rightarrow{}^o\!\!\!\cM'\circ{}^o\!\!\cL$,
and likewise for $\beta^{\dagger o}$. Say that
$\cL=(L_*,L^*,\eta^L,\eps^\cL)$, $\cL'=(L'_*,L'^*,\eta^{\cL'},\eps^{\cL'})$,
$\cM=(M_*,M^*,\eta^M,\eps^\cM)$ and $\cM'=(M'_*,M'^*,\eta^{\cM'},\eps^{\cM'})$.
Notice that for horizontally composable oriented squares of links we
have the identities :
\set\begin{equation}\label{eq_compose-horizont-sqlinks}
(\cD'\boxminus\cD)^\dagger=\cD'^\dagger\boxvert\cD^\dagger
\qquad
{}^\dagger(\cD'\boxminus\cD)={}^\dagger\!\cD\boxvert{}^\dagger\!\cD'.
\end{equation}
Thus, the operation $\cD\mapsto\cD^\dagger$ for oriented squares
of links differs in this respect from the analogous operation
for oriented squares in the $2$-category $\cA$ : see
\eqref{eq_change-orientation}.

With this notation, we attach to the oriented square $\cD$ the
{\em base change oriented square} in $\cA$:
$$
\Upsilon(\cD)\ :\quad
{\spreaddiagramcolumns{20pt}\diagram B' \ar[r]^-{M_*} \ar[d]_{L'^*}
\drtwocell\omit{_\ \ \ \ \ \Upsilon(\beta)} & B \ar[d]^{L^*} \\
A' \ar[r]_{M'_*} & A
\enddiagram}\qquad\qquad
$$
whose orientation is the {\em base change transformation}
$$
\Upsilon(\beta):=(\eps^\cL*M'_*L'^*)\odot(L^**\beta*L'^*)
\odot(L^*M_**\eta^{\cL'}):L^*M_*\Rightarrow M'_*L'^*.
$$

\begin{proposition}\label{prop_opp-links-and-base-ch}
With the notation of \eqref{subsec_base-change-map}, we have :
$$
\Upsilon(\cD)^o=\Upsilon(\cD^\dagger)
\qquad
{}^o\Upsilon(\cD)=\Upsilon({}^\dagger\!\cD).
$$
\end{proposition}
\begin{proof} The sought identity for $\Upsilon({}^\dagger\!\cD)$
will be proven by inspecting the following diagram :
$$
\xymatrix@C+40pt{B' \ar[r]^-{M_*} \ddouble
\drtwocell\omit{_\ \ \eps^\cM} & \drtwocell\omit{_\ \ \ \ \one_{M^*}}
B \ar[d]^{M^*} \rdouble & \drtwocell\omit{_\ \ \ \ \one_{M^*}}
B \ar[d]^{M^*} \rdouble & \drtwocell\omit{_\ \ \ \eta^\cM}
B \ar[d]^{M^*} \rdouble & \drtwocell\omit{_\ \ \ i_B}
B \ddouble \rdouble & B \ddouble \\
B' \rdouble \ar[d]_{L'^*} \drtwocell\omit{_\ \ \ \ \ \one_{L'^*}} &
B' \rdouble \ar[d]^{L'^*} \drtwocell\omit{_\ \ \ \ \ \one_{L'^*}} &
B' \rdouble  \ar[d]^{L'^*} \drtwocell\omit{_\ \ \ \ \eta^{\cL'}} &
B' \ar[r]|-{M_*} \ddouble \drtwocell\omit{_\ \ \ \ \one_{M_*}} &
B \rdouble \ddouble \drtwocell\omit{_\ \ \ i_B} & B \ddouble \\
A' \rdouble \ddouble \drtwocell\omit{_\ \ \ i_{A'}} &
A' \rdouble \ddouble \drtwocell\omit{_\ \ \ \ \one_{M'_*}} &
A' \ar[r]|-{L'_*} \ar[d]^{M'_*} \drtwocell\omit{_\ \ \beta} &
B' \ar[r]|-{M_*} \ar[d]^{M_*} \drtwocell\omit{_\ \ \ \ \one_{M_*}} &
B \rdouble \ddouble \drtwocell\omit{_\ \ \ i_B} & B \ddouble \\
A' \rdouble \ddouble \drtwocell\omit{_\ \ \ i_{A'}} &
A' \ar[r]|-{M'_*} \ddouble \drtwocell\omit{_\ \ \ \ \one_{M_*}} &
A \ar[r]|-{L_*} \ddouble \drtwocell\omit{_\ \ \ \ \eps^\cL} &
B \rdouble \ar[d]^{L^*} \drtwocell\omit{_\ \ \ \ \one_{L^*}} &
B \rdouble \ar[d]^{L^*} \drtwocell\omit{_\ \ \ \ \one_{L^*}} &
B \ar[d]^{L^*} \\
A' \rdouble \ddouble \drtwocell\omit{_\ \ \ i_{A'}} &
A' \ar[r]|-{M'_*} \ddouble \drtwocell\omit{_\ \ \ \ \eps^{\cM'}} &
A \rdouble \ar[d]^{M'^*} \drtwocell\omit{_\ \ \ \ \ \one_{M'^*}} &
A \rdouble  \ar[d]^{M'^*} \drtwocell\omit{_\ \ \ \ \ \one_{M'^*}} &
A \rdouble \ar[d]^{M'^*} \drtwocell\omit{_\ \ \ \ \eta^{\cM'}} &
A \ddouble \\
A' \rdouble & A' \rdouble & A' \rdouble & A' \rdouble &
A' \ar[r]_-{M'_*} & A.
}$$
Indeed, since the unit $\eta^{{}^o\!\!\cM}$ of ${}^o\!\!\cM$
is ${}^o(\eps^\cM)$ and the counit $\eps^{{}^o\!\!\cM'}$ of
${}^o\!\!\cM'$ is ${}^o(\eta^{\cM'})$, we have :
$$
\begin{aligned}
\Upsilon({}^\dagger\!\cD)&\,=(\eps^{{}^o\!\!\cM'}*{}^oL^*\ {}^o\!M_*)\odot
({}^o\!M'_**{}^o\!\beta^\dagger*{}^o\!M_*)\odot
({}^o\!M'_*\ {}^o\!L'^**\eta^{{}^o\!\cM}) \\
&\,={}^o((M'_*L'^**\eps^\cM)\odot(M'_**\beta^\dagger*M_*)
\odot(\eta^{\cM'}*L^*M_*))
\end{aligned}
$$
and if $\eta^{\cM\circ\cL'}$ is the unit of $\cM\circ\cL'$ and
$\eps^{\cL\circ\cM'}$ is the counit of $\cL\circ\cM'$, we have
$$
\begin{aligned}
\eta^{\cM\circ\cL'}=(M_**\eta^{\cL'}*M^*)\odot\eta^\cM
\qquad
\eps^{\cL\circ\cM'}=\eps^{\cM'}\odot(M'^**\eps^\cL*M'_*) \\
\beta^\dagger=(\eps^{\cL\circ\cM'}*L'^*M^*)\odot
(M'^*L^**\beta*L'^*M^*)\odot(M'^*L^**\eta^{\cM\circ\cL'}).
\end{aligned}
$$
Thus, if in the foregoing diagram we initially disregard the
first and last columns, we find that the composition of the
six remaining squares in the uppermost two rows equals
precisely $\eta^{\cM\circ\cL'}$, and the composition of the
six squares in the bottom two rows equals $\eps^{\cL\circ\cM'}$.
Therefore, if we further compose these two blocks of six squares
with the central row of three squares, we get precisely
$\beta^\dagger$. And if we finally compose the result with
the two omitted left and right columns, we then get
${}^o\Upsilon({}^\dagger\!\cD)$. But according to proposition
\ref{prop_square-algebra} we may compose these squares
in a different order : notice then that the composition
of the five squares in the uppermost row equals $\one_{M_*}$,
and the composition of the five square in the bottom row
equals $\one_{M'_*}$. So, we may disregard these two rows;
but in the remaining three central rows, we may likewise
disregard everything except the central column, since
the squares outside the central column are all identities.
So, finally we are left with the vertical composition of
the squares in the central column, and that gives precisely
$\Upsilon(\cD)$, as required. Lastly, we have :
$$
\begin{aligned}
\Upsilon(\cD^\dagger)&\,=(\eps^{\cM^o}*L'^{*o}\ M'^o_*)\odot
(M^o_**\beta^{\dagger o}*M'^o_*)\odot(M^o_*\ L^{*o}*\eta^{\cM'^o}) \\
&\,=((M'_*L'^**\eps^\cM)\odot(M'_**\beta^\dagger*M_*)
\odot(\eta^{\cM'}*L^*M_*))^o
\end{aligned}
$$
so the sought identity for $\Upsilon(\cD^\dagger)$ follows
from that for $\Upsilon({}^\dagger\!\cD)$.
\end{proof}

\begin{remark}\label{rem_when-Upsilon-inverts}
In the situation of \eqref{subsec_base-change-map}, suppose
that $\beta$ is an invertible $2$-cell of $\cA$ and that both
$L_*$ and $L'_*$ (resp. both $M_*$ and $M'_*$) are equivalences
in $\cA$. Then it follows easily from the definition (resp. and
from proposition \ref{prop_opp-links-and-base-ch}) that also
$\Upsilon(\beta)$ is an invertible $2$-cell.
\end{remark}

\begin{proposition}\label{prop_composition-of-sqlinks}
Consider two oriented squares of links :
$$
\cD\ :\quad
{\diagram A' \ar[r]^-{\cL'} \ar[d]_{\cM'}
\drtwocell\omit{_\ \beta} & B' \ar[d]^\cM \\
A \ar[r]_-\cL & B
\enddiagram}
\qquad\qquad\qquad
\cD'\ :\quad
{\diagram C' \ar[r]^-{\cP'} \ar[d]_{\cQ'}
\drtwocell\omit{_\ \alpha} & D' \ar[d]^\cQ \\
C \ar[r]_-\cP & D.
\enddiagram}$$
\begin{enumerate}
\item
Suppose that $A=C'$, $B=D'$ and $\cL=\cP'$, so that
$\cD'\boxvert\cD$ is well defined. Then :
$$
\Upsilon(\cD'\boxvert\cD)=\Upsilon(\cD')\boxminus\Upsilon(\cD).
$$
\item
Suppose that $B=C$, $B'=C'$ and $\cM=\cQ'$, so that
$\cD'\boxminus\cD$ is well defined. Then :
$$
\Upsilon(\cD'\boxminus\cD)=\Upsilon(\cD)\boxvert\Upsilon(\cD').
$$
\end{enumerate}
\end{proposition}
\begin{proof}(i): Say that $\cL=(L_*,L^*,\eta^\cL,\eps^\cL)$,
$\cP=(P_*,P^*,\eta^\cP,\eps^\cP)$, $\cQ=(Q_*,Q^*,\eta^\cQ,\eps^\cQ)$,
$\cM'=(M'_*,M'^*,\eta^{\cM'},\eps^{\cM'})$, $\cM=(M_*,M^*,\eta^\cM,\eps^\cM)$
and $\cQ'=(Q'_*,Q'^*,\eta^{\cQ'},\eps^{\cQ'})$. We consider the
following diagram :
$$
\xymatrix@C+40pt{B' \ar[r]^-{L'^*} \ddouble
\drtwocell\omit{^\eta^{\cL'}\ \ } & \drtwocell\omit{^\beta\ }
A' \ar[d]^{L'_*} \ar[r]^-{M'_*} & \drtwocell\omit{^\eps^\cL\ \ }
A \ar[d]^{L_*} \rdouble & \drtwocell\omit{^\one_{Q'^*}\ \ \ }
A \ar[r]^{Q'^*} \ddouble & \drtwocell\omit{^i_C\ \ }
C \ddouble \rdouble & C \ddouble \\
B' \rdouble \ddouble \drtwocell\omit{^i_{B'}\ \ } &
B' \ar[r]|-{M_*} \ddouble \drtwocell\omit{^\one_{M_*}\ \ \ } &
B \ar[r]|-{L^*} \ddouble \drtwocell\omit{^\eta^\cL\ \ } &
A \ar[r]|-{Q'^*} \ar[d]_{L_*} \drtwocell\omit{^\alpha\ } &
C \rdouble \ar[d]_{P_*} \drtwocell\omit{^\eps^\cP\ \ } & C \ddouble \\
B' \rdouble & B' \ar[r]_-{M_*} & B \rdouble & B \ar[r]_-{Q_*} &
D \ar[r]_-{P^*} & C
}$$
and notice that the composition of the five squares of the
top (resp. bottom) row equals $\Upsilon(\cD)$ (resp.
$\Upsilon(\cD')$), so the composition of all the square in
the diagram is precisely $\Upsilon(\cD')\boxminus\Upsilon(\cD)$.
On the other hand, the composition of the two squares in the
first (resp. second, resp. third, resp. fourth, resp. fifth)
column from the left equals $\eta^{\cL'}$ (resp. $\beta$, resp.
$\one_{L_*}$, resp. $\alpha$, resp. $\eps^\cP$), and then the
composition of these five columns also equals
$\Upsilon(\cD'\boxvert\cD)$, as stated.

(ii): We compute :
$$
\begin{aligned}
\Upsilon(\cD'\boxminus\cD)=
\Upsilon(\cD'^{\dagger\dagger}\boxminus\cD^{\dagger\dagger})
=\Upsilon((\cD'^\dagger\boxvert\cD^\dagger)^\dagger)
&\,=\Upsilon(\cD'^\dagger\boxvert\cD^\dagger)^o \\
&\,=(\Upsilon(\cD'^\dagger)\boxminus\Upsilon(\cD^\dagger))^o \\
&\,=\Upsilon(\cD^\dagger)^o\boxminus\Upsilon(\cD'^\dagger)^o \\
&\,=\Upsilon(\cD^{\dagger\dagger})\boxminus\Upsilon(\cD'^{\dagger\dagger}) \\
&\,=\Upsilon(\cD)\boxminus\Upsilon(\cD')
\end{aligned}
$$
where the second equality follows from
\eqref{eq_compose-horizont-sqlinks}, the third and sixth follow
from proposition \ref{prop_opp-links-and-base-ch}, and the fourth
follows from (i), and the fifth follows from
\eqref{eq_change-orientation}.
\end{proof}

\sset\subsubsection{}\label{subsec_murs-a-peche}
In the situation of \eqref{subsec_base-change-map}, the
links $\cL$ and $\cL'$ induce adjunctions $\theta^\cL$ and
$\theta^{\cL'}$ for the pairs of functors $((L^*)_*,(L_*)_*)$
and respectively $((L'^*)_*,(L'_*)_*)$, as in
\eqref{subsec_adjunction-from-link}. The base change
transformation $\Upsilon(\beta)$ relates these adjunctions,
as explained by the following :

\begin{proposition}\label{prop_relate-adjunctions-by-bc}
With the notation of \eqref{subsec_murs-a-peche}, for every
$(X\xrightarrow{f}A')\in\Ob(\cA/A')$, every
$(Y\xrightarrow{g}B')\in\Ob(\cA/B')$, and every $1$-cell
$(h,\alpha):g\to L'_*\circ f$ of $\cA/B'$ we have the identity :
$$
(M'_**\theta^{\cL'}_{f,g}(h,\alpha))\circ(\one_Y,\Upsilon(\beta)*g)=
\theta^\cL_{M'_*\circ f,M_*\circ g}((\one_X,\beta*f)\circ(h,M_**\alpha))
\qquad
\text{in $\cA/A$}.
$$
\end{proposition}
\begin{proof} The left-hand side of the stated identity is
computed as the composition of the following oriented squares :
$$
\xymatrix@C+25pt{ Y \ar[r]^-h \ar[d]_g \ddrtwocell\omit{^\alpha\ } &
X \rdouble \ar[d]_f \drtwocell\omit{^\one\ } &
X \ar[d]^f \rrdouble \ddrrtwocell\omit{^\one\ } & & X \dddouble \\
B' \ddouble & A' \rdouble \ar[d]_{L'_*} \drtwocell\omit{^\eps^{\cL'}\ } &
A' \ddouble & & \\
B' \rdouble \ddouble \drtwocell\omit{^\one\ } &
B' \ar[r]^-{L'^*} \ddouble \drtwocell\omit{^\eta^{\cL'}\ } &
A' \ar[r]^-{M'_*} \ar[d]_{L'_*} \drtwocell\omit{^\beta\ } &
A \ar[d]_{L_*} \rdouble \drtwocell\omit{^\eps^\cL\ } & A \ddouble \\
B' \rdouble & B' \rdouble & B' \ar[r]^-{M_*} & B \ar[r]^-{L^*} & A.
}$$
By applying as usual proposition \ref{prop_square-algebra},
and the triangular identities for the pair $(\eta^{\cL'},\eps^{\cL'})$,
we easily see that the same composition of squares yields as
well the right-hand side of the sought identity : the details
are left to the reader.
\end{proof}

\sset\subsubsection{}\label{subsec_define-pseudo-funct-Upsilon}
Recall that the links of $\cA$ are the objects of the $2$-category
$2\tdu\sMorph(\sLink(\cA))$ : see example \ref{ex_2-arrows}(i). An
oriented square $\cD$ as in \eqref{subsec_base-change-map}
yields a $1$-cell in this $2$-category
$$
\cD(\beta,\cM,\cM'):\cL'\to\cL.
$$
Given a pair of oriented squares
$\cD_1(\beta_1,\cM_1,\cM'_1),\cD_2(\beta_2,\cM_2,\cM'_2):
\cL'\to\cL$, a $2$-cell $(\alpha,\alpha'):\cD_1\Rightarrow\cD_2$
is the datum of a pair of $2$-cells of the $2$-category
$\sLink(\cA)$
\set\begin{equation}\label{eq_more-space}
\alpha:\cM_1\Rightarrow\cM_2
\quad
\alpha':\cM'_1\Rightarrow\cM'_2
\quad\text{such that}\quad
\beta_2\odot(\alpha*\cL')=(\cL*\alpha')\odot\beta_1.
\end{equation}

\begin{theorem}\label{th_Upsilon-is-pseudo-fctr}
The rules :
$$
\cL:=(L_*,L^*,\eta,\eps)\mapsto L^{*o}
\qquad
\cD\mapsto\Upsilon(\cD)^o
\qquad
(\alpha,\alpha')\mapsto(\alpha^o,\alpha'^o)
$$
for every link $\cL$, every oriented square $\cD$ of
$\sLink(\cA)$, and every $2$-cell $(\alpha,\alpha')$ as
in \eqref{subsec_define-pseudo-funct-Upsilon} define a
strict pseudo-functor :
$$
\Upsilon_\cA:2\tdu\sMorph(\sLink(\cA))^o\to 2\tdu\sMorph(\cA^o).
$$
\end{theorem}
\begin{proof} The functoriality of $\Upsilon$ for composition
of $1$-cells is already clear from proposition
\ref{prop_composition-of-sqlinks}. Let us next check that for
every $2$-cell $(\alpha,\alpha'):\cD_1\Rightarrow\cD_2$ as in
\eqref{subsec_define-pseudo-funct-Upsilon}, the pair
$(\alpha^o,\alpha'^o)$ is a $2$-cell
$\Upsilon(\cD_1)^o\Rightarrow\Upsilon(\cD_2)^o$
in $2\tdu\sMorph(\cA^o)$. We come down to showing the identity
\set\begin{equation}\label{eq_kidneys}
\Upsilon(\beta_2)^o\odot(\alpha^o*L'^{*o})=
(L'^{*o}*\alpha')\odot\Upsilon(\beta_1)^o.
\end{equation}
To this aim, for every $X\in\Ob(\cA)$ let $\bone_X:X\to X$
be the link $(\one_X,\one_X,i_X,i_X)$ (where $i_X$ is the
identity automorphism of $\one_X$); also, for $i=1,2$ let
$\cP_i:=(P_{i*}:X\to Y,P^*_i,\eta_i,\eps_i)$ be a link of
$\cA$, and $\tau:P_{1*}\Rightarrow P_{2*}$ a $2$-cell
of $\cA$. We consider the following oriented squares of
$\sLink(\cA)$ and of $\cA^o$ :
$$
\cD_\tau\ :\
{\diagram X \ar[r]^-{\bone_X} \ar[d]_{\cP_2}
\drtwocell\omit{_\ \tau} & X \ar[d]^{\cP_1} \\
Y \ar[r]_-{\bone_Y} & Y
\enddiagram}
\qquad\qquad
\cE_\tau\ :\
{\diagram Y^o \ar[r]^-{\one_{Y^o}} \ar[d]_{P^o_{2*}}
\drtwocell\omit{_\ \tau^o} &
Y^o \ar[d]^{P^o_{1*}} \\
X^o \ar[r]_-{\one_{X^o}} & X^o.
\enddiagram}$$
With this notation, by applying proposition
\ref{prop_opp-links-and-base-ch} it is easily seen that
$$
\Upsilon(\cD_\tau)^o=\cE_\tau
$$
and \eqref{eq_kidneys} is equivalent to the identity :
$$
\Upsilon(\cD_1)^o\boxminus\cE_{\alpha'}=
\cE_\alpha\boxminus\Upsilon(\cD_2)^o.
$$
In view of \eqref{eq_compose-horizont-sqlinks} and proposition
\ref{prop_composition-of-sqlinks} we are then reduced to checking
that
$$
\cD_1\boxminus\cD_{\alpha'}=\cD_\alpha\boxminus\cD_2
$$
which follows immediately from \eqref{eq_more-space}. Lastly,
the functoriality of $\Upsilon$ for the two types of compositions
of $2$-cells is obvious from the definition.
\end{proof}

\sset\subsubsection{}\label{subsec_variant-of-Upsilon}
Resume the situation of \eqref{subsec_base-change-map}, and
notice that the units and counits for $\cM$ and $\cM'$ do not
play any role in the definition of $\Upsilon(\beta)$, and
neither do the $1$-cells $M^*$ and $M'^*$. This leads to the
following variant of theorem \ref{th_Upsilon-is-pseudo-fctr}.
We replace $2\tdu\sMorph(\sLink(\cA))$ by the $2$-category
$$
\wLink(\cA)
$$
of {\em weak links in $\cA$}, whose objects are still the links of
$\cA$; for any two links $\cL:=(L_*:A\to B,L^*,\eta^\cL,\eps^\cL)$
and $\cL':=(L'_*:A'\to B',L'^*,\eta^{\cL'},\eps^{\cL'})$, the
$1$-cells $\cL\to\cL'$ are the $1$-cells
$(M_*,M'_*,\beta):L_*\to L'_*$ of the $2$-category
$2\tdu\sMorph(\cA)$, {\em i.e.} the oriented squares
$$
\cD\ :\quad
{\diagram A' \ar[r]^-{L'_*} \ar[d]_{M'_*}
\drtwocell\omit{_\ \beta} & B' \ar[d]^{M_*} \\
A \ar[r]_-{L_*} & B.
\enddiagram}$$
The $2$-cells of $\wLink(\cA)$, and the composition laws for
$1$-cells and $2$-cells  are then the same as for
$2\tdu\sMorph(\cA)$. For any such $1$-cell $\cD$, we may then
define the $2$-cell $\Upsilon(\beta)$ and the oriented square
$\Upsilon(\cD)$ as in \eqref{subsec_base-change-map}. A direct
inspection shows that part (i) of proposition
\ref{prop_composition-of-sqlinks} still holds for weak links,
with the same proof. The proof of part (ii) of proposition
\ref{prop_composition-of-sqlinks} does not apply to weak links,
but we may check the assertion directly. Indeed, consider another
$1$-cell $\cP':=(P'_*,P'^*,\eta^{\cP'},\eps^{\cP'})
\to\cP:=(P_*,P^*,\eta^\cP,\eps^\cP)$ of $\wLink(\cA)$ :
$$
\cD'\ :\quad
{\diagram B' \ar[r]^-{M_*} \ar[d]_{Q_*}
\drtwocell\omit{_\ \alpha} & C' \ar[d]^{Q_*} \\
B \ar[r]_-{P_*} & C.
\enddiagram}$$
Set :
$$\begin{aligned}
X:=&\,\eps^\cL*M'_*L'^*P'^* & & & Y:=&\,L^*P^*Q_**\eta^{\cP'} \\
\delta:=&\,(P_**\beta)\odot(\alpha*L'_*) &
\quad\alpha'&\,:=P^**\alpha*P'^*\quad & \beta':=&\,L^**\beta*L'^*.
\end{aligned}
$$
By definition we have :
$$
\begin{aligned}
&\Upsilon(\cD'\boxminus\cD)=X\!\odot\!
(L^*\!*\!\eps^\cP\!*\!L_*M'_*L'^*P'^*)
\!\odot\!(L^*P^*\!*\!\delta\!*\!L'^*P'^*)\!\odot\!
(L^*P^*Q_*P'_*\!*\!\eta^{\cL'}\!*\!P'^*)
\!\odot\!Y \\
&\Upsilon(\cD)\boxvert\Upsilon(\cD')=X\odot
((\beta'\odot(L^*M_**\eta^{\cL'}))*P'^*)\odot
(L^**((\eps^\cP*M_*P'^*)\odot\alpha'))\odot Y.
\end{aligned}
$$
Hence we are reduced to checking that
$$
Z:=(\eps^\cP*L_*M'_*L'^*)\odot(P^*\!*\!\delta\!*\!L'^*)\odot
(P^*Q_*P'_*\!*\!\eta^{\cL'})
$$
equals
$$
Z':=((\beta*L'^*)\odot(M_**\eta^{\cL'}))\odot(\eps^\cP*M_*)
\odot(P^**\alpha).
$$
The latter follows by applying repeatedly remark
\ref{rem_equiv-2-cat}(i) : the details shall be left to the reader.
With these preliminaries, we may then also repeat the proof of
theorem \ref{th_Upsilon-is-pseudo-fctr} : summing up, we obtain
a strict pseudo-functor :
$$
\Upsilon:\wLink(\cA)^o\to 2\tdu\sMorph(\cA^o)
$$
given by the same rules as for the previously defined one on
$2\tdu\sMorph(\sLink(\cA))^o$.

\begin{remark}\label{rem_double-upsilon}
In the situation of \eqref{subsec_base-change-map}, notice that
the base change oriented square $\Upsilon(\cD)$ can be regarded
as a $1$-cell $\cM\to\cM'$ in the $2$-category of weak links.
By the discussion of \eqref{subsec_variant-of-Upsilon}, the
oriented square $\Upsilon(\Upsilon(\cD))$ is then well defined. A
straightforward computation that we shall leave to the reader yields
the identity :
$$
\Upsilon(\Upsilon(\cD))=\beta^\dagger.
$$
\end{remark}

\sset\subsubsection{}\label{subsec_transfer-base-change}
Let us consider two $1$-cells in the $2$-category of weak
links in $\cA$ :
$$
\begin{aligned}
(M'_{1*},M_{1*},\delta_1):\cL'_1:=(L'_{1*},L'^*_1,\eta^{\cL'_1},\eps^{\cL'_1})\to
\cL_1:=(L_{1*},L^*_1,\eta^{\cL_1},\eps^{\cL_1}) \\
(M'_{2*},M_{2*},\delta_2):\cL'_2:=(L'_{2*},L'^*_2,\eta^{\cL'_2},\eps^{\cL'_2})\to
\cL_2:=(L_{2*},L^*_2,\eta^{\cL_2},\eps^{\cL_2})
\end{aligned}
$$
and a diagram of $2$-cells :
$$
\xymatrix{ A'_2 \ar[rrrrr]^-{L'_{2*}} \ar[ddd]_{M'_{2*}} & & &
\dltwocell\omit{_\ \alpha} & &
B'_2 \ar[ddd]^{M_{2*}} \\ 
& \dltwocell\omit{^\ \phi'} & A'_1 \ar[r]_-{L'_{1*}} \ar[d]_{M'_{1*}}
\ar[llu]_{Q'_{1*}} \drtwocell\omit{_\ \delta_1} &
B'_1 \ar[d]^{M_{1*}} \ar[rru]^{Q'_{2*}} &
\drtwocell\omit{_\ \phi} \\
& & A_1 \ar[r]^-{L_{1*}} \ar[lld]_{Q_{1*}} \drtwocell\omit{_\ \beta} &
B_1 \ar[rrd]^{Q_{2*}} & & \\
A_2 \ar[rrrrr]_-{L_{2*}} & & & & & B_2
}$$
so that $\alpha$ and $\beta$ are $1$-cells $\cL'_1\to\cL'_2$
and $\cL_1\to\cL_2$ in $\wLink(\cA)$. We complete the diagram
by adding as well the $2$-cell $\delta_2$ (the reader should
picture this as a cubical diagram whose faces are oriented
by the given $2$-cells), and we suppose moreover that the
resulting diagram {\em commutes on $2$-cells, i.e.} we have :
$$
(\delta_2\boxvert\alpha)\boxminus\phi'=
\phi\boxminus(\beta\boxvert\delta_1).
$$
We may then regard $\alpha$ and $\beta$ as $1$-cells of
$\wLink(\cA)$ :
$$
(Q'_{1*},Q'_{2*},\alpha):\cL'_1\to\cL'_2
\qquad
(Q_{1*},Q_{2*},\beta):\cL_1\to\cL_2
$$
and the pair $(\phi',\phi)$ as a $2$-cell of $\wLink(\cA)$ :
$$
(\phi',\phi):(M'_{2*},M_{2*},\delta_2)\circ(Q'_{1*},Q'_{2*},\alpha)
\Rightarrow(Q_{1*},Q_{2*},\beta)\circ(M'_{1*},M_{1*},\delta_1).
$$
The pseudo-functor $\Upsilon$ of \eqref{subsec_variant-of-Upsilon}
then yields the $2$-cell :
$$
(\phi'^o,\phi^o):\Upsilon(M'_{2*},M_{2*},\delta_2)^o\circ
\Upsilon(Q'_{1*},Q'_{2*},\alpha)^o\Rightarrow
\Upsilon(Q_{1*},Q_{2*},\beta)^o\circ\Upsilon(M'_{1*},M_{1*},\delta_1)^o
$$
which in turns translates as the diagram of $2$-cells
$$
\xymatrix@C+15pt{ B'_2 \ar[rrrrr]^-{M_{2*}} \ar[ddd]_{L'^*_2} & & &
\dltwocell\omit{^\ \phi} & &
B_2 \ar[ddd]^{L^*_2} \\ 
& \dltwocell\omit{^\ \ \ \ \Upsilon(\alpha)} &
B'_1 \ar[r]_-{M_{1*}} \ar[d]_{L'^*_1}
\ar[llu]_{Q'_{2*}} \drtwocell\omit{_\ \ \ \ \ \ \Upsilon(\delta_1)} &
B_1 \ar[d]^{L^*_1} \ar[rru]^{Q_{2*}} &
\drtwocell\omit{_\ \ \ \ \ \Upsilon(\beta)} \\
& & A'_1 \ar[r]^-{M'_{1*}} \ar[lld]_{Q'_{1*}} \drtwocell\omit{^\phi'} &
A_1 \ar[rrd]^{Q_{1*}} & & \\
A'_2 \ar[rrrrr]_-{M'_{2*}} & & & & & A_2
}$$
completed to a cubical diagram, by adding in the
$2$-cell $\Upsilon(M'_{2*},M_{2*},\delta_2)$, and commuting again
on $2$-cells, {\em i.e.} verifying the identity :
$$
(\Upsilon(\beta)\boxminus\Upsilon(\delta_1))\boxvert\phi=
\phi'\boxvert(\Upsilon(\delta_2)\boxminus\Upsilon(\alpha)).
$$

\begin{remark}\label{rem_transit-base-change}
(i)\ \
In the situation of \eqref{subsec_transfer-base-change}, suppose
that $\Upsilon(\alpha)$, $\Upsilon(\beta)$, $\Upsilon(\delta_1)$,
$\phi$ and $\phi'$ are invertible $2$-cells of $\cA$. In this case,
we deduce that the same holds for $\Upsilon(\delta_2)*Q'_2$, and if
$Q'_2$ is an equivalence in $\cA$, we conclude easily that also
$\Upsilon(\delta_2)$ is an invertible $2$-cell.

(ii)\ \
For instance, consider the case where $L_{1*}=L_{2*}$, $L'_{1*}=L'_{2*}$,
$M_{1*}=M_{2*}$, $M'_{1*}=M'_{2*}$ and where all the $1$-cells
$Q_{1*},Q_{2*},Q'_{1*},Q'_{2*}$ and all the $2$-cells $\alpha$,
$\beta$, $\phi$ and $\phi'$ are identities (and then
$\delta_1=\delta_2$). In this case, applying (i) we conclude
that $\Upsilon(\delta_1:\cL'_1\to\cL_1)$ is invertible if and
only if the same holds for $\Upsilon(\delta_2:\cL'_2\to\cL_2)$.

(iii)\ \
We may also consider the variant where we invert the
orientation of the $2$-cells $\phi$ and $\phi'$ (and
keep unchanged the others); then the commutativity on
$2$-cells becomes the condition :
$$
\phi\boxminus(\delta_2\boxvert\alpha)=
(\beta\boxvert\delta_1)\boxminus\phi'
$$
so that the pair $(\phi,\phi')$ can be regarded as a
$1$-cell of $\wLink(\cA)$ as in the foregoing, but with
reversed direction. After applying $\Upsilon$, we deduce
the identity
$$
\phi'\boxvert(\Upsilon(\beta)\boxminus\Upsilon(\delta_1))=
(\Upsilon(\delta_2)\boxminus\Upsilon(\alpha))\boxvert\phi.
$$
Especially, (i) applies {\em verbatim} to this variant as well.
\end{remark}

\sset\subsubsection{}\label{subsec_functoriality-of-Link}
Lastly, let $\cA$ and $\cB$ be two $2$-categories, and
$\phi:\cA\to\cB$ a pseudo-functor with coherence constraints
$(\delta_\bullet,\gamma_\bullet)$. Then $\phi$ induces a
pseudo-functor :
$$
\sLink(\phi):\sLink(\cA)\to\sLink(\cB)
\qquad
A\mapsto \phi A
$$
that assigns to every link $\cL:=(F,G,\eta^\cL,\eps^\cL):A\to B$
of $\cA$ the datum
$$
\phi(\cL):=(\phi F,\phi G,\eta^{\phi(\cL)},\eps^{\phi(\cL)})
$$
where $\eta^{\phi(\cL)}$ and $\eps^{\phi(\cL)}$ are the unique
$2$-cells that make commute the diagrams
$$
\xymatrix{ \one_{\phi B} \ar@{=>}[rr]^-{\eta^{\phi(\cL)}}
\ar@{=>}[d]_{\delta_B} & &
\phi F\circ\phi G \ar@{=>}[d]^{\gamma_{F,G}} &
\phi G\circ\phi F \ar@{=>}[rr]^-{\eps^{\phi(\cL)}}
\ar@{=>}[d]_{\gamma_{G,F}} & &
\one_{\phi A} \ar@{=>}[d]^{\delta_A} \\
\phi\one_B \ar@{=>}[rr]^-{\phi(\eta^\cL)} & & \phi(F\circ G) &
\phi(G\circ F) \ar@{=>}[rr]^-{\phi(\eps^\cL)} & & \phi\one_A.
}$$
We leave to the reader the verification that $\phi(\cL)$ is
indeed a link $\phi A\to\phi B$ in $\cB$. To every transformation
of links $\beta:\cL\Rightarrow\cL'$, the pseudo-functor
$\sLink(\phi)$ assigns the transformation
$\phi(\beta):\phi(\cL)\Rightarrow\phi(\cL')$. The associativity
constraint of $\sLink(\phi)$ assigns to every composable pair
of morphisms of links $\cL:=(F,G,\eta^\cL,\eps^\cL):A\to B$,
$\cL':=(F',G',\eta^{\cL'},\eps^{\cL'}):B\to C$ the $2$-cell
$\gamma^*_{\cL,\cL'}:=\gamma_{F,F'}:
\phi(\cL')\circ\phi(\cL)\Rightarrow\phi(\cL'\circ\cL)$. Likewise,
the unit constraint of $\sLink(\phi)$ assigns to every object
$A$ of $\cA$ the $2$-cell $\delta^*_A:=\delta_A$.

\begin{proposition}\label{prop_functoriality-of-Link}
In the situation of \eqref{subsec_functoriality-of-Link}, the
following diagram commutes :
$$
\xymatrix{2\tdu\sMorph(\sLink(\cA))^o \ar[rr]^-{\Upsilon_\cA}
\ar[d]_{2\tdu\sMorph(\sLink(\phi))^o} & &
2\tdu\sMorph(\cA^o) \ar[d]^{2\tdu\sMorph(\phi^o)} \\
2\tdu\sMorph(\sLink(\cB))^o \ar[rr]^-{\Upsilon_\cB} & &
2\tdu\sMorph(\cB^o).
}$$
\end{proposition}
\begin{proof} The commutativity on objects and $2$-cells of
$2\tdu\sMorph(\sLink(\cA))^o$ is immediate from the definition.
To check commutativity on $1$-cells, consider a diagram $\cD$
as in \eqref{subsec_base-change-map}; we come down to verifying
that the $2$-cell of $\cB$ :
$$
X:=(\eps^{\phi(\cL)}*\phi M'_*\phi L'^*)\odot(\phi L^**
(\gamma^{-1}_{M'_*,L_*}\odot\phi(\beta)\odot\gamma_{L'_*,M_*})*\phi L'^*)
\odot(\phi L^*\phi M_**\eta^{\phi(\cL')})
$$
equals the $2$-cell :
$$
Y:=\gamma_{M'_*,L'^*}^{-1}\odot\phi((\eps^\cL*M'_*L'^*)\odot(L^**\beta*L'^*)
\odot(L^*M_**\eta^{\cL'}))\odot\gamma_{L^*,M_*}.
$$
To this aim, set
$$
\begin{aligned}
Z:=&\,((\delta^{-1}_A\odot\phi(\eps^\cL)\odot\gamma_{L^*,L_*})*\phi M'_*)
\odot(\phi L^**\gamma^{-1}_{M'_*,L_*}) \\
Z':=&\,(\gamma_{L'_*,M_*}*\phi L'^*)\odot
(\phi M_**(\gamma^{-1}_{L'_*,L'^*}\odot\phi(\eta^{\cL'})\odot\delta_{B'}))
\end{aligned}
$$
so that
$X=(Z*\phi L'^*)\odot X'\odot(\phi L^**Z')$, with
$X':=\phi L^**\phi(\beta)*\phi L'^*$. We compute :
$$
\begin{aligned}
Z=&\,\gamma_{\one_A,M'_*}\odot(\phi(\eps^\cL)*\phi M'_*)\odot
(\gamma_{L^*,L_*}*\phi M'_*)\odot(\phi L^**\gamma^{-1}_{M'_*,L_*}) \\
=&\,\gamma_{\one_A,M'_*}\odot(\phi(\eps^\cL)*\phi M'_*)\odot
\gamma_{L^*L_*,M'_*}^{-1}\odot\gamma_{L^*,L_*M'_*} \\
=&\phi(\eps^\cL*M'_*)\odot\gamma_{L^*,L_*M'_*} \\
Z'=&\,(\gamma_{L'_*,M_*}*\phi L'^*)\odot(\phi M_**\gamma^{-1}_{L'_*,L'^*})
\odot(\phi M_**\phi(\eta^{\cL'}))\odot\gamma^{-1}_{M_*,\one_B} \\
=&\,\gamma^{-1}_{M_*L'_*,L'^*}\odot\gamma_{M_*,L'_*L'^*}\odot
(\phi M_**\phi(\eta^{\cL'}))\odot\gamma^{-1}_{M_*,\one_B} \\
=&\gamma^{-1}_{M_*L'_*,L'^*}\odot\phi(M_**\eta^{\cL'}).
\end{aligned}
$$
On the other hand, we have :
$$
\begin{aligned}
\gamma_{M'_*,L'^*}^{-1}\odot\phi(\eps^\cL*M'_*L'^*)=&\,
(\phi(\eps^\cL*M'_*)*\phi L'^*)\odot\gamma^{-1}_{L^*L_*M'_*,L'^*} \\
\phi(L^*M_**\eta^{\cL'})\odot\gamma_{L^*,M_*}=&\,
\gamma_{L^*,M_*L'_*L'^*}\odot(\phi L^**\phi(M_**\eta^{\cL'}))
\end{aligned}
$$
so we are reduced to showing the identity :
$$
(\gamma_{L^*,L_*M'_*}*\phi L'^*)\odot X'
\odot(\phi L^**\gamma^{-1}_{M_*L'_*,L'^*})=\gamma^{-1}_{L^*L_*M'_*,L'^*}
\odot\phi(L^**\beta*L'^*)\odot\gamma_{L^*,M_*L'_*L'^*}.
$$
The latter is the same as the identity :
$\gamma_{L^*,L_*M'_*,L'^*}\odot X'=\phi(L^**\beta*L'^*)\odot
\gamma_{L^*,M_*L'_*,L'^*}$, which holds by virtue of remark
\ref{rem_pseudo-funct}(iii).
\end{proof}

\sset\subsubsection{}\label{subsec_ups-and-equiv}
We conclude with a further construction derived from the
formalism of base change, that shall help us in the following
section. First, to every $2$-category $\cA$ we attach the
$2$-category
$$
\Equiv(\cA)
$$
defined as the sub-$2$-category of $2\tdu\sMorph(\cA)$
whose objects are the equivalences in $\cA$ (see definition
\ref{def_adjoint-1-cells}(iii)); for any two equivalences
$f:X\to Y$ and $f':X'\to Y'$ in $\cA$, the $1$-cells $f\to f'$
in $\Equiv(\cA)$ are the essentially commutative oriented squares :
\set\begin{equation}\label{eq_this-is-it}
{\diagram X \ar[r]^-f \ar[d]_{g'} \drtwocell\omit{_\ \beta} &
Y \ar[d]^g \\
X' \ar[r]_-{f'} & Y'
\enddiagram}\end{equation}
and the $2$-cells between such $1$-cells are as in
$2\tdu\sMorph(\cA)$. We define a pseudo-functor
$$
\sL:\Equiv(\cA)\to\wLink(\cA)
$$
as follows. For any $(f:X\to Y)\in\Ob(\Equiv(\cA))$ we choose
a quasi-inverse $f^\dagger:Y\to X$ (see  definition
\ref{def_adjoint-1-cells}(iii)) and the unit $\eta^f$ and
counit $\eps^f$ of an adjunction for $(f^\dagger,f)$, and we set
$$
\sL(f):=(f^\dagger,f,\eta^f,\eps^f)\in\Ob(\wLink(\cA)).
$$
Then every $1$-cell $(g,g',\beta):f\to f'$ as in
\eqref{eq_this-is-it} yields a $1$-cell $\sL(f)\to\sL(f')$
in $\wLink(\cA)$ which we denote $\sL(g,g',\beta)$; likewise,
$\sL$ is the identity on $2$-cells. Then of course the coherence
constraints of $\sL$ are given by identities, but $\sL$
{\em is not} necessarily strict, since we do not necessarily
have $(g\circ f)^\dagger=f^\dagger\circ g^\dagger$ for every
$f,g\in\Ob(\Equiv(\cA))$.

\sset\subsubsection{}\label{subsec_tilde-upsilon}
Moreover, we have a natural strict isomorphism of
$2$-categories :
$$
\Equiv(\cA)\isom\Equiv(\cA^o)^o.
$$
Namely, to every equivalence $f:X\to Y$ of $\cA$ we assign
the equivalence $f^o:Y^o\to X^o$ of $\cA^o$, and to every
$1$-cell \eqref{eq_this-is-it} we assign the following $1$-cell
$f'^o\to f^o$ of $\Equiv(\cA^o)^o$ :
$$
\xymatrix@C+30pt{ Y'^o \ar[r]^-{f'^o} \ar[d]_{g^o}
\drtwocell\omit{_\ \ \ \ \ \ \ (\beta^o)^{-1}} & X'^o \ar[d]^{g'^o} \\
Y^o \ar[r]_-{f^o}  & X^o.
}$$
Lastly, let $(h,h',\beta'):f\to f'$ be another $1$-cell of
$\Equiv(\cA)$, and
$(\alpha_1,\alpha_2):(g,g',\beta)\Rightarrow(h,h',\beta')$
any $2$-cell; recall that the latter is a pair of $2$-cells
$\alpha_1:g'\Rightarrow h'$ and $\alpha_2:g\Rightarrow h$ of
$\cA$ such that $\alpha_1\odot\beta=\beta'\odot\alpha_2$.
Then our strict isomorphism assigns to $(\alpha_1,\alpha_2)$
the pair $(\alpha_2^o,\alpha_1^o)$, which is easily seen to
be a $2$-cell $(g'^o,g^o,(\beta^o)^{-1})\Rightarrow
(h'^o,h^o,(\beta'^o)^{-1})$ of $\Equiv(\cA^o)^o$.

Next, by composing the pseudo-functors $\sL$ of
\eqref{subsec_ups-and-equiv} and $\Upsilon$ of
\eqref{subsec_variant-of-Upsilon} we get a pseudo-functor
$\Upsilon^o\circ\sL:\Equiv(\cA)\to 2\tdu\sMorph(\cA^o)^o$,
and taking into account remark \ref{rem_when-Upsilon-inverts}
we see that $\Upsilon^o\circ\sL$ factors through the inclusion
strict pseudo-functor $\Equiv(\cA^o)^o\to 2\tdu\sMorph(\cA^o)^o$.
Therefore, we can further compose with the foregoing
isomorphism to deduce a pseudo-functor
$$
\tilde\Upsilon:\Equiv(\cA)\to\Equiv(\cA)
\qquad
f\mapsto f^\dagger
\qquad
(g,g',\beta)\mapsto\Upsilon(g,g',\beta)^{-1}.
$$

\subsection{Adjunctions in \texorpdfstring{$2$}{2}-categories}
\label{sec_2-adjoints}
We begin by introducing two constructions that will often allow
us to reduce the proof of assertions concerning general
pseudo-functors, to the special case of strict, or at least
{\em unital} pseudo-functors, in the sense of definition
\ref{def_unital-pseudo}.

\begin{definition}\label{def_unital-pseudo}
Let $\cA$ and $\cB$ be two $2$-categories, and $F:\cA\to\cB$
a pseudo-functor with coherence constraint
$(\delta^F_\bullet,\gamma^F_{\bullet\bullet})$. We say that $F$ is
{\em unital} if $F\one_A=\one_{FA}$ and $\delta^F_A=i_{FA}$ for
every $A\in\Ob(\cA)$. Since the component $\delta^F_\bullet$ is
determined by the these conditions, we shall usually say that
the coherence constraint of a unital pseudo-functor $F$ is
the component $\gamma^F_{\bullet\bullet}$.
\end{definition}

\begin{remark}\label{rem_unital}
Let $F:\cA\to\cB$ be a unital pseudo-functor with coherence
constraint $\gamma^F$.

(i)\ \ 
Directly from the composition axiom we see that for every
pair of $1$-cells $A\xrightarrow{\ f\ }B\xrightarrow{\ g\ }C$
of $\cA$, the $2$-cell $\gamma^F_{f,g}$ equals $\one_{Fg}$ if
$f=\one_A$, and equals $\one_{Ff}$ if $g=\one_B$.

(ii)\ \
Let $G:\cA\to\cB$ be another unital pseudo-functor, and
$\alpha:F\Rightarrow G$ any pseudo-natural transformation
with coherence constraint $\tau$. From the first coherence
axiom for $\alpha$ it is clear that $\tau_{\one_A}=\one_{\alpha_A}$
for every $A\in\Ob(\cA)$.
\end{remark}

\begin{proposition}\label{prop_towards-2-yoneda}
For every pseudo-functor $F:\cA\to\cB$ there exists a unital
pseudo-functor $F^u:\cA\to\cB$ with a pseudo-natural isomorphism
$F\isom F^u$.
\end{proposition}
\begin{proof} Let $(\delta^F,\gamma^F)$ be the coherence
constraint of $F$; for every $1$-cell $f:A\to B$ of $\cA$,
we define the $2$-cell of $\cB$
$$
\omega_f:=\left\{ \begin{array}{ll}
                  \one_{Ff} & \text{if $f\neq\one_A$} \\
                  \delta^F_A & \text{if $f=\one_A$}
                  \end{array}\right.
\qquad
F^uf\Rightarrow Ff.
$$
We shall show more precisely, that there exists a unique
well-defined unital pseudo-functor $F^u:\cA\to\cB$ with
coherence constraint $\gamma^{F^u}$ such that :
\begin{itemize}
\item
$F^uA=FA$ for every $A\in\Ob(\cA)$
\item
$F^uf=Ff$ for every $A,B\in\Ob(\cA)$ every $1$-cell $f:A\to B$
such that $f\neq\one_A$
\item
for every pair of $1$-cells
$A\xrightarrow{\ f\ }B\xrightarrow{\ g\ }C$ of $\cA$ we have
$$
\gamma^{F^u}_{f,g}=
\omega^{-1}_{g\circ f}\odot\gamma^F_{f,g}\odot(\omega_g*\omega_f)
$$
\item
for every pair of $1$-cells $f,f':A\to B$ of $\cA$ and every
$2$-cell $\beta:f\Rightarrow f'$, we have
$$
F^u\beta=\omega_{f'}^{-1}\odot F\beta\odot\omega_f.
$$
\end{itemize}
Indeed, the uniqueness of such $F^u$ is obvious; it is also
clear that for every $A,B\in\Ob(\cA)$ the stated rules yield
a well-defined functor $F^u_{AB}:\cA(A,B)\to\cB(F^uA,F^uB)$.
For the naturality of $\gamma^{F^u}$ consider $1$-cells
$f,f':A\to B$ and $g,g':B\to C$ in $\cA$ and $2$-cells
$\beta_1:f\Rightarrow f'$ and $\beta_2:g\Rightarrow g'$;
we notice that
$$
\begin{aligned}
F^u\beta_2*F^u\beta_1=&\,
(\omega_{g'}^{-1}\odot F\beta_2\odot\omega_g)*
(\omega_{f'}^{-1}\odot F\beta_1\odot\omega_f) \\
=&\,
(\omega_{g'}^{-1}*\omega_{f'}^{-1})\odot
(F\beta_2*F\beta_1)\odot(\omega_g*\omega_f)
\end{aligned}
$$
whence
$$
\begin{aligned}
F^u(\beta_2*\beta_1)\odot\gamma^{F^u}_{f,g}=&\,
\omega^{-1}_{g'\circ f'}\odot F(\beta_2*\beta_1)\odot
\gamma^F_{f,g}\odot(\omega_g*\omega_f) \\
=&\,
\omega^{-1}_{g'\circ f'}\odot\gamma^F_{f',g'}\odot
(F\beta_2*F\beta_1)\odot(\omega_g*\omega_f) \\
=&\,
\gamma^{F^u}_{f',g'}\odot(F^u\beta_2*F^u\beta_1)
\end{aligned}
$$
which is the contention, by remark \ref{rem_pseudo-funct}(ii).
Next, let us check the composition axiom for $\gamma^{F^u}$;
we consider three $1$-cells
$A\xrightarrow{\ f\ }B\xrightarrow{\ g\ }C\xrightarrow{\ h\ }D$
of $\cA$, and we need to show the identity
$$
X:=\gamma^{F^u}_{f,h\circ g}\odot(\gamma^{F^u}_{g,h}*F^uf)=
Y:=\gamma^{F^u}_{g\circ f,h}\odot(F^uh*\gamma^{F^u}_{f,g}).
$$
We compute, on the one hand :
$$
\begin{aligned}
X=&\,
\omega^{-1}_{h\circ g\circ f}\odot\gamma^F_{f,h\circ g}\odot
(\omega_{h\circ g}\odot\omega_f)\odot
((\omega^{-1}_{h\circ g}\odot\gamma^F_{g,h}*(\omega_h*\omega_g))*F^uf) \\
=&\,
\omega^{-1}_{h\circ g\circ f}\odot\gamma^F_{f,h\circ g}\odot
(\omega_{h\circ g}*\omega_f)\odot(\omega^{-1}_{h\circ g}*F^uf)
\odot((\gamma^F_{g,h}\odot(\omega_h*\omega_g))*F^uf) \\
=&\,
\omega^{-1}_{h\circ g\circ f}\odot\gamma^F_{f,h\circ g}\odot
(F(h\circ g)*\omega_f)\odot
((\gamma^F_{g,h}\odot(\omega_h*\omega_g))*F^uf) \\
=&\,
\omega^{-1}_{h\circ g\circ f}\odot\gamma^F_{f,h\circ g}\odot
(F(h\circ g)*\omega_f)\odot(\gamma^F_{g,h}*F^uf)\odot
(\omega_h*\omega_g*F^uf) \\
=&\,
\omega^{-1}_{h\circ g\circ f}\odot\gamma^F_{f,h\circ g}\odot
(\gamma^F_{g,h}*Ff)\odot(Fh*Fg*\omega_f)\odot
(\omega_h*\omega_g*F^uf) \\
=&\,
\omega^{-1}_{h\circ g\circ f}\odot\gamma^F_{f,h\circ g}\odot
(\gamma^F_{g,h}*Ff)\odot(\omega_h*\omega_g*\omega_f) \\
=&\,
\omega^{-1}_{h\circ g\circ f}\odot\gamma^F_{g\circ f,h}\odot
(Fh*\gamma^F_{f,g})\odot(\omega_h*\omega_g*\omega_f)
\end{aligned}
$$
and on the other hand :
$$
\begin{aligned}
Y=&\,
\omega^{-1}_{h\circ g\circ f}\odot\gamma^F_{g\circ f,h}\odot
(\omega_h*\omega_{g\circ f})\odot(F^uh*(\omega^{-1}_{g\circ f}
\odot\gamma^F_{f,g}\odot(\omega_g*\omega_f))) \\
=&\,
\omega^{-1}_{h\circ g\circ f}\odot\gamma^F_{g\circ f,h}\odot
(\omega_h*\omega_{g\circ f})\odot(F^uh*\omega^{-1}_{g\circ f})
\odot(F^uh*(\gamma^F_{f,g}\odot(\omega_g*\omega_f))) \\
=&\,
\omega^{-1}_{h\circ g\circ f}\odot\gamma^F_{g\circ f,h}\odot
(\omega_h*F(g\circ f))\odot
(F^uh*\gamma^F_{f,g})\odot(F^uh*\omega_g*\omega_f) \\
=&\,
\omega^{-1}_{h\circ g\circ f}\odot\gamma^F_{g\circ f,h}\odot
(Fh*\gamma^F_{f,g})\odot(\omega_h*Fg*Fh)\odot
(F^uh*\omega_g*\omega_f) \\
=&\,
\omega^{-1}_{h\circ g\circ f}\odot\gamma^F_{g\circ f,h}\odot
(Fh*\gamma^F_{f,g})\odot(\omega_h*\omega_g*\omega_f)
\end{aligned}
$$
whence the contention. Lastly, the unit axiom can be
verified by a simple inspection.

The sought pseudo-isomorphism is given by the rule that
assigns $\alpha_A:=\one_{FA}:FA\to F^uA$ to every $A\in\Ob(\cA)$,
and the oriented square
$$
{\diagram FA \ar[r]^-{\alpha_A} \ar[d]_{Ff}
\drtwocell\omit{_\ \ \omega_f} & F^uA \ar[d]^{F^uf} \\
FB \ar[r]_-{\alpha_B} & F^uB
\enddiagram}
\qquad
\text{to every $1$-cell $f:A\to B$ of $\cA$}.
$$
The naturality of $\omega_f$ and the coherence axioms follow
by a simple inspection.
\end{proof}

\begin{remark}\label{rem_uniPsFun}
(i)\ \
For any two $2$-categories $\cA$ and $\cB$, let us denote by
$$
\uniPsFun(\cA,\cB)
$$
the sub-$2$-category of $\sPsFun(\cA,\cB)$ whose objects
are the unital pseudo-functors $\cA\to\cB$, and whose
$\Hom$-categories are given by the categories $\sPsNat(F,G)$,
for any two such unital pseudo-functors $F,G$. For every
pseudo-functor $F:\cA\to\cB$, let also $\alpha^F:F\isom F^u$
be the pseudo-natural isomorphism furnished by proposition
\ref{prop_towards-2-yoneda}. It is easily seen that the rules
$F\mapsto F^u$ for every pseudo-functor $F:\cA\to\cB$,
$(\beta:F\Rightarrow G)\mapsto
\beta^u:=\alpha^G\odot\beta\odot(\alpha^F)^{-1}$
for every pseudo-natural transformation $\beta$, and
$(\Xi:\beta\leadsto\beta')\mapsto\Xi^u:=\alpha^G*\Xi*(\alpha^F)^{-1}$
for every modification $\Xi$, define a strict $2$-equivalence
of $2$-categories
$$
(-)^u:\sPsFun(\cA,\cB)\isom\uniPsFun(\cA,\cB)
$$
(details left to the reader). Likewise, if
$i:\uniPsFun(\cA,\cB)\to\sPsFun(\cA,\cB)$ denotes the inclusion
strict pseudo-functor, the rule : $F\mapsto\alpha^F$ for every
pseudo-functor $F:\cA\to\cB$ yields a strict pseudo-natural
isomorphism of strict pseudo-functors :
\set\begin{equation}\label{eq_really-strict}
\one_{\sPsFun(\cA,\cB)}\isom i\circ(-)^u.
\end{equation}

(ii)\ \
Let $F:\cA\to\cB$ and $G:\cB\to\cC$ be two
pseudo-functors, $\alpha^F:F\isom F^u$, $\alpha^G:G\isom G^u$,
$\alpha^{G\circ F}:G\circ F\isom(G\circ F)^u$ the corresponding
pseudo-natural isomorphisms as in (i). Since $\alpha^F$ and
$\alpha^G$ are not strict, we have distinct pseudo-natural
isomorphisms
$$
(\alpha^G*F^u)\odot(G*\alpha^F),(G^u*\alpha^F)\odot(\alpha^G*F):
G\circ F\isom G^u\circ F^u
$$
that are related by an invertible modification (see
example \ref{ex_modifications}(ii)). By composing with
$(\alpha^{G\circ F})^{-1}$, we get therefore two pseudo-natural
isomorphisms of unital pseudo-functors
$$
(G\circ F)^u\isom G^u\circ F^u
$$
but in general these pseudo-functors are different.
\end{remark}

\begin{proposition}\label{prop_strictification}
For every pseudo-functor $F:\cA\to\cB$ there exist :
\begin{enumerate}
\alphaenu
\item
A $2$-category $\cA^F$ and strict pseudo-functors
$\cA\xleftarrow{\pi^F}\cA^F\xrightarrow{F^\flat}\cB$.
\item
With an isomorphism of pseudo-functors
$\omega^F:F\circ\pi^F\isom F^\flat$.
\item
A pseudo-functor $\sigma^F:\cA\to\cA^F$ such that :
$$
\pi^F\circ\sigma^F=\one_\cA
\qquad
F^\flat\circ\sigma^F=F
\qquad
\omega^F*\sigma^F=\one_F.
$$
\item
A pseudo-natural equivalence $\mu^F:\sigma^F\circ\pi^F\isom\one_{\cA^F}$.
\end{enumerate}
\end{proposition}
\begin{proof} We let $\cA^F$ be the $2$-category with
$\Ob(\cA^F)=\Ob(\cA)$, and whose $1$-cells
$$
(f,g,\lambda):X\to Y
\qquad
\text{for every $X,Y\in\Ob(\cA)$}
$$
are all the data consisting of a $1$-cell $f:X\to Y$ of $\cA$,
a $1$-cell $g:FX\to FY$ of $\cB$, and an invertible $2$-cell
$\lambda:g\isom Ff$. Given two such $1$-cells
$(f,g,\lambda),(f',g',\lambda'):X\to Y$, the $2$-cells
$$
(f,g,\lambda)\Rightarrow(f',g',\lambda')
$$
are the pairs $(\alpha,\beta)$ where $\alpha:f\Rightarrow f'$
(resp. $\beta:g\Rightarrow g'$) is a $2$-cell of $\cA$ (resp.
of $\cB$), such that
$$
\lambda'\odot\beta=F\alpha\odot\lambda.
$$
The composition of two $1$-cells $X\xrightarrow{(f_1,g_1,\lambda_1)}
Y\xrightarrow{(f_2,g_2,\lambda_2)}Z$ is the $1$-cell
$$
(f_2,g_2,\lambda_2)\circ(f_1,g_1,\lambda_1):=
(f_2\circ f_1,g_2\circ g_1,\gamma^F_{f_1,f_2}\odot(\lambda_2*\lambda_1))
$$
where $(\delta^F,\gamma^F)$ denotes the coherence constraint
of the pseudo-functor $F$. The associativity of this composition
law follows easily from the composition axiom for $\gamma^F$ : the
details are left to the reader. Especially, for every
$X\in\Ob(\cA)$ the triple $(\one_X,\one_{FX},\delta^F_X)$
is the identity $1$-cell of $X$ in $\cA^F$. Given three
$1$-cells $(f,g,\lambda),(f',g',\lambda'),(f'',g'',\lambda''):
X\to Y$, and two $2$-cells
$(\alpha,\beta):(f,g,\lambda)\Rightarrow(f',g',\lambda')$ and
$(\alpha',\beta'):(f',g',\lambda')\Rightarrow(f'',g'',\lambda'')$,
we let
$$
(\alpha',\beta')\odot(\alpha,\beta):=
(\alpha'\odot\alpha,\beta'\odot\beta)
$$
which is a well defined $2$-cell
$(f,g,\lambda)\Rightarrow(f'',g'',\lambda'')$. Lastly, if
$(f_1,g_1,\lambda_1),(f'_1,g'_1,\lambda'_1):X\to Y$ and
$(f_2,g_2,\lambda_2),(f'_2,g'_2,\lambda'_2):Y\to Z$ are four
given $1$-cells, and $(\alpha_1,\beta_1):(f_1,g_1,\lambda_1)
\Rightarrow(f'_1,g'_1,\lambda'_1)$,
$(\alpha_2,\beta_2):(f_2,g_2,\lambda_2)\Rightarrow
(f'_2,g'_2,\lambda'_2)$ are two $2$-cells, we set :
$$
(\alpha_2,\beta_2)*(\alpha_1,\beta_1):=
(\alpha_2*\alpha_1,\beta_2*\beta_1)
$$
which is a well defined $2$-cell
$(f_2,g_2,\lambda_2)\circ(f_1,g_1,\lambda_1)\Rightarrow
(f'_2,g'_2,\lambda'_2)\circ(f'_1,g'_1,\lambda'_1)$. The
associativity of these two composition laws are obvious, and
it is then easily seen that $\cA^F$ is a well defined
$2$-category with such laws for $1$-cells and $2$-cells : the
details are left to the reader.

The sought strict pseudo-functors $F^\flat$ and $\pi^F$ are
given by the rules :
$$
\begin{aligned}
F^\flat X&\,:=FX
& & &
\pi^FX&\,:=X
& & &
\text{for every $X\in\Ob(\cA^F)$} \\
F^\flat(f,g,\lambda)&\,:=g
& & &
\pi^F(f,g,\lambda)&\,:=f
& & &
\text{for every $1$-cell $(f,g,\lambda)$} \\
F^\flat(\alpha,\beta)&\,:=\beta
& & &
\pi^F(\alpha,\beta)&\,:=\alpha
& & &
\text{for every $2$-cell $(\alpha,\beta)$}.
\end{aligned}
$$
The isomorphism $\omega^F$ is given by the rule : $X\mapsto\one_{FX}$
for every $X\in\Ob(\cA^F)$, and its coherence costraint assigns
to every $1$-cell $(f,g,\lambda):X\to Y$ the oriented square
$$
\xymatrix{ FX \rdouble \ar[d]_{Ff}
\drtwocell\omit{_\lambda} & FX \ar[d]^g \\
FY \rdouble & FY.
}$$

Next, the pseudo-functor $\sigma^F$ is given by the rules :
$$
X\mapsto X
\qquad
f\mapsto(f,Ff,\one_{Ff})
\qquad
\alpha\mapsto(\alpha,F\alpha)
$$
for every $X\in\Ob(\cA)$, every $1$-cell $f$, and every $2$-cell
$\alpha$ of $\cA$. The coherence constraints of $\sigma^F$ are
given by the rules :
$$
X\mapsto(i_X,\delta^F_X)
\qquad
(f,f')\mapsto(\one_{f'\circ f},\gamma^F_{f,f'})
$$
for every $X\in\Ob(\cA)$ and every composable pair of $1$-cells
$(f,f')$ of $\cA$. The required coherence axioms are easily
verified, and the desired identities as in (c) follow
straightforwardly : details left to the reader. Lastly, we
define $\mu^F_X:=\one_{\sigma^FX}=(\one_X,\one_{FX},\delta^F_X)$
for every $X\in\Ob(\cA^F)$. For every $1$-cell
$(f,g,\lambda):X\to Y$ of $\cA^F$, notice that
$\sigma^F\circ\pi^F(f,g,\lambda)=(f,Ff,\one_{Ff})$; then the
coherence constraint of $\mu^F$ assigns to $(f,g,\lambda)$
the oriented square in $\cA^F$ :
$$
\xymatrix@C+40pt{ X \ar[r]^-{(\one_X,\one_{FX},\delta^F_X)}
\ar[d]_{(f,Ff,\one_{Ff})} \drtwocell\omit{_\ \ \ \ \ \ \ (\one_f,\lambda)}
& X \ar[d]^{(f,g,\lambda)} \\
Y \ar[r]_-{(\one_Y,\one_{FY},\delta^F_Y)} & Y.
}$$
A little diagram chase that we leave to the reader shows
that these rules yield the sought pseudo-natural equivalence.
\end{proof}

\begin{proposition}\label{prop_strictify-pairs}
{\em (i)}\ \
Let $\cA,\cB$ be two $2$-categories, $F:\cA\to\cB$ and
$G:\cB\to\cA$ two pseudo-functors. Then there exist :
\begin{enumerate}
\alphaenu
\item
Two $2$-categories $\cA^{F,G}$ and $\cB^{G,F}$, and a diagram
of strict pseudo-functors :
$$
\xymatrix@C+20pt{
\cA & \cA^{F,G} \ar@<.5ex>[r]^-{F^\flat} \ar[l]_-{\pi^{F,G}} &
\cB^{G,F} \ar@<.5ex>[l]^-{G^\flat} \ar[r]^-{\pi^{G,F}} & \cB.
}$$
\item
With two isomorphisms of pseudo-functors :
$$
\omega^{F,G}:F\circ\pi^{F,G}\isom\pi^{G,F}\circ F^\flat
\qquad\text{and}\qquad
\omega^{G,F}:G\circ\pi^{G,F}\isom\pi^{F,G}\circ G^\flat.
$$
\item
Two pseudo-functors $\sigma^{F,G}:\cA\to\cA^{F,G}$
and $\sigma^{G,F}:\cB\to\cB^{G,F}$ such that :
$$
\begin{aligned}
\pi^{F,G}\circ\sigma^{F,G}&\,=\one_\cA
& \qquad\qquad &
F^\flat\circ\sigma^{F,G}&\,=\sigma^{G,F}\circ F
& \qquad\qquad &
\omega^{F,G}*\sigma^{F,G}&\,=\one_F \\
\pi^{G,F}\circ\sigma^{G,F}&\,=\one_\cB
& \qquad\qquad &
G^\flat\circ\sigma^{G,F}&\,=\sigma^{F,G}\circ G
& \qquad\qquad &
\omega^{G,F}*\sigma^{G,F}&\,=\one_G.
\end{aligned}
$$
\item
With two isomorphisms of pseudo-functors :
$$
\psi^{F,G}:\sigma^{F,G}\circ\pi^{F,G}\isom\one_{\cA^{F,G}}
\quad\text{and}\quad
\psi^{G,F}:\sigma^{G,F}\circ\pi^{G,F}\isom\one_{\cB^{G,F}}
$$
\end{enumerate}

{\em (ii)}\ \
Suppose moreover that $F$ and $G$ are unital. Then we have
as well :
$$
G^\flat*\psi^{G,F}=(\psi^{F,G}*G^\flat)\odot(\sigma^{F,G}*\omega^{G,F})
\qquad
F^\flat*\psi^{F,G}=(\psi^{G,F}*F^\flat)\odot(\sigma^{G,F}*\omega^{F,G}).
$$
\end{proposition}
\begin{proof}(i): To begin with, for every $k,n\in\N$ let us set :
$$
H_{2k}:=G
\qquad
H_{2k+1}:=F
\qquad
K_0:=\one_\cA
\qquad
K_{n+1}:=H_{n+1}\circ H_{n-1}\circ\cdots\circ H_1.
$$
Then we let $\cA^{F,G}$ be the $2$-category with
$\Ob(\cA^{F,G})=\Ob(\cA)$, and whose $1$-cells
$$
(f_\bullet,\lambda_\bullet):X\to Y
\qquad
\text{for every $X,Y\in\Ob(\cA)$}
$$
are all the systems of $1$-cells $(f_n:K_nX\to K_nY~|~n\in\N)$,
and of invertible $2$-cells
$(\lambda_n:f_{n+1}\isom H_{n+1}f_n~|~n\in\N)$. Given such
$1$-cells
$(f_\bullet,\lambda_\bullet),(f'_\bullet,\lambda'_\bullet):X\to Y$,
the $2$-cells
$$
\alpha_\bullet:(f_\bullet,\lambda_\bullet)\Rightarrow
(f'_\bullet,\lambda'_\bullet)
$$
are the systems of $2$-cells $(\alpha_n:f_n\to f'_n~|~k\in\N)$
such that :
\set\begin{equation}\label{eq_recursive}
H_{n+1}(\alpha_n)\odot\lambda_n=\lambda'_n\odot\alpha_{n+1}
\qquad
\text{for every $n\in\N$}.
\end{equation}
The composition of two $1$-cells
$X\xrightarrow{(f_{1,\bullet},\lambda_{1,\bullet})}
Y\xrightarrow{(f_{2,\bullet},\lambda_{2,\bullet})}Z$ is the $1$-cell
$$
(f_{2,\bullet},\lambda_{1,\bullet})\circ(f_{1,\bullet},\lambda_{1,\bullet}):=
(f_{2,n}\circ f_{1,n},\gamma^{H_{n+1}}_{f_{1,n},f_{2,n}}\odot
(\lambda_{2,n}*\lambda_{1,n})~|~n\in\N)
$$
where $(\delta^{H_n},\gamma^{H_n})$ denotes the coherence
constraint of the pseudo-functor $H_n$, for every $n\in\N$.
Just as in the proof of proposition \ref{prop_strictification},
the associativity of this rule follows easily from the composition
axioms for $\gamma^{H_n}$, and for every $X\in\Ob(\cA)$,
the identity $1$-cell of $X$ in $\cA^{F,G}$ is the datum
$\one_X^{F,G}:=(\one_{K_nX},\delta^{H_{n+1}}_{K_nX}~|~n\in\N)$.
Given three $1$-cells
$(f_\bullet,\lambda_\bullet),(f'_\bullet,\lambda'_\bullet),
(f''_\bullet,\lambda''_\bullet):X\to Y$, and two $2$-cells
$\alpha_\bullet:(f_\bullet,\lambda_\bullet)\Rightarrow
(f'_\bullet,\lambda'_\bullet)$ and $\alpha'_\bullet:
(f'_\bullet,\lambda'_\bullet)\Rightarrow(f''_\bullet,\lambda''_\bullet)$,
we let
$$
\alpha'_\bullet\odot\alpha_\bullet:=
(\alpha'_n\odot\alpha_n~|~n\in\N)
$$
which is a well defined $2$-cell
$(f_\bullet,\lambda_\bullet)\Rightarrow(f''_\bullet,\lambda''_\bullet)$.
Lastly, if
$(f_{1,\bullet},\lambda_{1,\bullet}),
(f'_{1,\bullet},\lambda'_{1,\bullet}):X\to Y$ and
$(f_{2,\bullet},\lambda_{2,\bullet}),
(f'_{2,\bullet},\lambda'_{2,\bullet}):Y\to Z$ are four
given $1$-cells, and
$\alpha_{1,\bullet}:(f_{1,\bullet},\lambda_{1,\bullet})
\Rightarrow(f'_{1,\bullet},\lambda'_{1,\bullet})$,
$\alpha_{2,\bullet}:(f_{2,\bullet},\lambda_{2,\bullet})\Rightarrow
(f'_{2,\bullet},\lambda'_{2,\bullet})$ are two $2$-cells, we set :
$$
\alpha_{2,\bullet}*\alpha_{1,\bullet}:=
(\alpha_{2,n}*\alpha_{1,n}~|~n\in\N)
$$
which is a well defined $2$-cell
$(f_{2,\bullet},\lambda_{2,\bullet})\circ
(f_{1,\bullet},\lambda_{1,\bullet})\Rightarrow
(f'_{2,\bullet},\lambda'_{2,\bullet})\circ
(f'_{1,\bullet},\lambda'_{1,\bullet})$. The associativity of these
two composition laws are obvious, and it is then easily seen
that $\cA^{F,G}$ is a well defined $2$-category with such laws for
$1$-cells and $2$-cells : the details are left to the reader.

We define $\cB^{G,F}$ by exactly the same rules, after swapping
$\cA$ and $F$ with respectively $\cB$ and $G$. Then the sought
strict pseudo-functors $F^\flat$ and $\pi^{F,G}$ are
given by the rules :
$$
\begin{aligned}
F^\flat X&\,:=FX
& & &
\pi^{F,G}X&\,:=X
& & &
\text{for every $X\in\Ob(\cA^{F,G})$} \\
F^\flat(f_\bullet,\lambda_\bullet)&\,:=(f_{n+1},\lambda_{n+1}~|~n\in\N)
& & &
\pi^{F,G}(f_\bullet,\lambda_\bullet)&\,:=f_0
& & &
\text{for every $1$-cell $(f_\bullet,\lambda_\bullet)$} \\
F^\flat(\alpha_\bullet)&\,:=(\alpha_{n+1}~|~n\in\N)
& & &
\pi^{F,G}(\alpha_\bullet)&\,:=\alpha_0
& & &
\text{for every $2$-cell $\alpha_\bullet$}
\end{aligned}
$$
and again, $G^\flat$ and $\pi^{G,F}$ are defined by the same
rules, after swapping $F$ and $G$.

The isomorphism $\omega^{F,G}$ is given by the rule :
$X\mapsto\one_{FX}$ for every $X\in\Ob(\cA^F)$, and its
coherence costraint assigns to every $1$-cell
$(f_\bullet,\lambda_\bullet):X\to Y$ the oriented square
$$
\xymatrix{ FX \rdouble \ar[d]_{Ff_0}
\drtwocell\omit{_\lambda_0} & FX \ar[d]^{f_1} \\
FY \rdouble & FY.
}$$
Lastly, the pseudo-functor $\sigma^{F,G}$ is given by the rules :
$$
X\mapsto X
\qquad
f\mapsto(K_nf,\one_{K_{n+1}f}~|~n\in\N)
\qquad
\alpha\mapsto(K_n\alpha~|~n\in\N)
$$
for every $X\in\Ob(\cA)$, every $1$-cell $f$, and every $2$-cell
$\alpha$ of $\cA$.
For every $n\in\N$, let $(\delta^{K_n},\gamma^{K_n})$ be the
coherence constraint of $K_n$; the coherence constraint of
$\sigma^{F,G}$ is given by the rules :
$$
X\mapsto(\delta^{K_n}_X~|~n\in\N)
\qquad
(f,f')\mapsto(\gamma^{K_n}_{f,f'}~|~n\in\N)
$$
for every $X\in\Ob(\cA)$ and every composable pair of $1$-cells
$(f,f')$ of $\cA$, where $H$. The required coherence axioms are
verified as in the proof of the corresponding assertions of
the proposition \ref{prop_strictification}, and shall again
be left to the reader. Likewise, one defines $\omega^{G,F}$
and $\sigma^{G,F}$ by the same rules, after swapping the roles
of $F$ and $G$. Then the desired identities in (c) follow
straightforwardly. It remains to exhibit an isomorphism of
pseudo-functors $\psi^{F,G}:\sigma^{F,G}\circ\pi^{F,G}\isom\one_{\cA^{F,G}}$.
To this aim, notice that $\sigma^{F,G}\circ\pi^{F,G}$ is given
by the rules :
$$
X\mapsto X
\qquad
(f_\bullet,\lambda_\bullet)\mapsto(K_nf_0,\one_{K_{n+1}f_0}~|~n\in\N)
\qquad
\alpha_\bullet\mapsto(K_n\alpha_0~|~n\in\N)
$$
for every object $X$, every $1$-cell $(f_\bullet,\lambda_\bullet)$
and every $2$-cell $\alpha_\bullet$ of $\cA^{F,G}$; then we let
$\psi^{F,G}_X:=\one^{F,G}_X$ for every such $X$, and to every such
$(f_\bullet,\lambda_\bullet)$ we attach the coherence constraint
$$
\tau^{\psi^{F,G}}_{(f_\bullet,\lambda_\bullet),\bullet}:(f_\bullet,\lambda_\bullet)
\to(K_nf_0,\one_{K_{n+1}f_0}~|~n\in\N)
$$
defined inductively by the rules :
$$
\tau^{\psi^{F,G}}_{(f_\bullet,\lambda_\bullet),0}:=\one_{f_0}
\qquad
\tau^{\psi^{F,G}}_{(f_\bullet,\lambda_\bullet),n+1}:=
H_{n+1}(\tau^{\psi^{F,G}}_{(f_\bullet,\lambda_\bullet),n})\odot\lambda_n
\qquad
\text{for every $n\in\N$}.
$$
The coherence axioms for $\tau^\psi$ amount to the identities :
$$
\begin{aligned}
\tau^{\psi^{F,G}}_{\one^{F,G}_X,n}&\,=\delta^{K_n}_X \\
A_n:=\gamma^{K_n}_{f_0,f'_0}\odot
(\tau^{\psi^{F,G}}_{(f'_\bullet,\lambda'_\bullet),n}*
\tau^{\psi^{F,G}}_{(f_\bullet,\lambda_\bullet),n})&\,=
B_n:=\tau^{\psi^{F,G}}_{(f''_\bullet,\lambda''_\bullet),n}
\end{aligned}
\qquad
\text{for every $n\in\N$}
$$
for every object $X$ of $\cA$ and every composable pair of $1$-cells
$(f_\bullet,\lambda_\bullet)$ and $(f'_\bullet,\lambda'_\bullet)$, with
$(f''_\bullet,\lambda''_\bullet):=(f'_\bullet,\lambda'_\bullet)\circ
(f_\bullet,\lambda_\bullet)$. The first identity is checked easily
by induction on $n$, recalling that
$\delta^{K_{n+1}}_X=H_{n+1}(\delta^{K_n}_X)\odot\delta^{H_{n+1}}_{K_nX}$
for every $n\in\N$. The second identity is obvious for $n=0$.
Suppose then that the second identity is known for some $n\in\N$;
we compute :
$$
\begin{aligned}
A_{n+1}&\,=
H_{n+1}(\gamma^{K_n}_{f_0,f'_0})\odot\gamma^{H_{n+1}}_{K_nf_0,K_nf'_0}\odot
((H_{n+1}(\tau^\psi_{(f'_\bullet,\lambda'_\bullet),n})\odot\lambda'_n)*
(H_{n+1}(\tau^\psi_{(f_\bullet,\lambda_\bullet),n})\odot\lambda_n)) \\
&\,=H_{n+1}(\gamma^{K_n}_{f_0,f'_0})\odot\gamma^{H_{n+1}}_{K_nf_0,K_nf'_0}\odot
(H_{n+1}(\tau^\psi_{(f'_\bullet,\lambda'_\bullet),n})*
H_{n+1}(\tau^\psi_{(f_\bullet,\lambda_\bullet),n}))\odot(\lambda'_n*\lambda_n) \\
&\,=H_{n+1}(\gamma^{K_n}_{f_0,f'_0})\odot
H_{n+1}(\tau^\psi_{(f'_\bullet,\lambda'_\bullet),n}*
\tau^\psi_{(f_\bullet,\lambda_\bullet),n})\odot\gamma^{H_{n+1}}_{f_0,f'_0}
\odot(\lambda'_n*\lambda_n) \\
&\,=H_{n+1}(A_n)\odot\lambda''_n \\
&\,=H_{n+1}(B_n)\odot\lambda''_n \\
&\,=B_{n+1}
\end{aligned}
$$
where the first equality follows after recalling that
$\gamma^{K_{n+1}}_{f_0,f'_0}=H_{n+1}(\gamma^{K_n}_{f_0,f'_0})\odot
\gamma^{H_{n+1}}_{K_nf_0,K_nf'_0}$. The same rules yield, {\em mutatis
mutandis}, the sought isomorphism
$\psi^{G,F}:\sigma^{G,F}\circ\pi^{G,F}\isom\one_{\cB^{G,F}}$.

(ii): Let us set as well :
$$
H'_{2k}:=F
\qquad
H'_{2k+1}:=G
\qquad
K'_0:=\one_\cB
\qquad
K'_{n+1}:=H'_{n+1}\circ H'_{n-1}\circ\cdots\circ H'_1.
$$
With this notation, since $F$ are unital, we have :
$$
\begin{aligned}
(G^\flat*\psi^{G,F})_X&\,=G^\flat(\one^{G,F}_X)
=(\one_{K'_{n+1}X},i_{K'_{n+2}X}~|~n\in\N) \\
(\psi^{F,G}*G^\flat)_X&\,=\psi^{F,G}_{GX}
=(\one_{K_nGX},i_{K_{n+1}GX}~|~n\in\N) \\
(\sigma^{F,G}*\omega^{G,F})_X&\,=\sigma^{F,G}_{GX}
=(\one_{K_nGX},i_{K_{n+1}GX}~|~n\in\N).
\end{aligned}
$$
Noticing that $K'_{n+2}=K_{n+1}G$ and
$\gamma^{H_{n+1}}_{\one_{K_nGX},\one_{K_nGX}}=\delta^{H_{n+1}}_{K_nGX}=
i_{K_nGX}$, we deduce :
$$
(\psi^{F,G}*G^\flat)_X\!\circ\!(\sigma^{F,G}\!*\!\omega^{G,F})_X
=(G^\flat*\psi^{G,F})_X
\qquad
\text{for every $X\in\Ob(\cB^{G,F})$}.
$$
Next, the coherence constraint of $\sigma^{F,G}*\omega^{G,F}$
is given by the rule :
$$
(X\xrightarrow{(f_\bullet,\lambda_\bullet)}Y)\mapsto
((\gamma^{K_n}_{Gf_0,\one_{{GY}}})^{-1}\odot K_n\lambda_0\odot
\gamma^{K_n}_{\one_{GX},f_1}~|~n\in\!\N)=(K_n\lambda_0~|~n\in\!\N)
$$
whereas the ones of $\psi^{F,G}*G^\flat$ and respectively
$G^\flat*\psi^{G,F}$ are defined by the rules :
$$
(f_\bullet,\lambda_\bullet)\mapsto
(\tau^{\psi^{F,G}}_{(f_{n+1},\lambda_{n+1}~|~n\in\N),n}~|~n\in\N)
\qquad
(f_\bullet,\lambda_\bullet)\mapsto
(\tau^{\psi^{G,F}}_{(f_\bullet,\lambda_\bullet),n+1}~|~n\in\N).
$$
Hence we are reduced to showing that
$$
K_n\lambda_0\odot\tau^{\psi^{F,G}}_{(f_{n+1},\lambda_{n+1}~|~n\in\N),n}=
\tau^{\psi^{G,F}}_{(f_\bullet,\lambda_\bullet),n+1}
\qquad
\text{for every $n\in\N$}.
$$
The latter follows by an easy induction on $n$ : details left
to the reader.
\end{proof}

\begin{definition}\label{def_using-mods}
Let $\cA$, $\cB$ be two $2$-categories, and $F,F':\cA\to\cB$
two pseudo-functors.

(i)\ \
We say that a pseudo-natural transformation $F\Rightarrow F'$
is a {\em pseudo-natural equivalence\/} if it is an equivalence
in the $2$-category $\sPsFun(\cA,\cB)$ (in the sense of definition
\ref{def_adjoint-1-cells}(iii)).

(ii)\ \
We say that $F$ is {\em fully faithful} (resp. {\em strongly
faithful}) if for every $A,A'\in\Ob(\cA)$, the functor
$F_{AA'}:\cA(A,A')\to\cB(FA,FA')$ is an equivalence (resp. is
an isomorphism of categories).

(iii)\ \
We say that $\cA$ is a {\em full sub-$2$-category} (resp. a
{\em strong sub-$2$-category}) of $\cB$, if it is a sub-$2$-category
and the inclusion pseudo-functor $\cA\to\cB$ is fully (resp.
strongly) faithful.

(iv)\ \
We say that $F$ is a {\em $2$-equivalence\/} (resp. a
{\em strong $2$-equivalence}) from $\cA$ to $\cB$ if it is
fully faithful (resp. strongly faithful), and for every
$B\in\Ob(\cB)$ there exists $A\in\Ob(\cA)$ with an equivalence
(resp. an isomorphism) $FA\to B$.
\end{definition}

\begin{remark}\label{rem_pseudo-commutative}
Consider a square diagram of $2$-categories and pseudo-functors
$$
\cD\qquad : \qquad
{\diagram \cA \ar[r]^-F \ar[d]_H & \cB \ar[d]^G \\
\cC \ar[r]^K & \cD.
\enddiagram}
$$
We shall say that $\cD$ is {\em essentially commutative}
(resp. {\em pseudo-commutative}) if there exists an isomorphism
of pseudo-functors (resp. a pseudo-natural equivalence)
$\alpha:G\circ F\Rightarrow K\circ H$. Even though the
set of (small) $2$-categories does not form a $2$-category,
the datum $(\cD,\alpha)$ is a type of oriented square analogous
to those contemplated in \eqref{subsec_square-algebra}. Especially,
two such data $(\cD,\alpha)$ and $(\cD',\alpha')$ can be composed
if they share a side, in exactly the same way as for usual
oriented square. Then, it is clear that a horizontal or
vertical composition of essentially commutative (resp.
pseudo-commutative) squares is again essentially commutative
(resp. pseudo-commutative). On the other hand, the discussion
of remark \ref{rem_pseudo-natural}(iv,v) shows that the
identities of proposition \ref{prop_square-algebra} hold
only up to isomorphism of pseudo-natural transformations,
{\em i.e.} up to invertible modifications.
\end{remark}

\begin{lemma}\label{lem_equiv-is-preserved}
Let $\cA,\cB$ be two $2$-categories, $F:\cA\to\cB$ a
pseudo-functor, and $f:A\to B$ a $1$-cell in $\cA$.
We have :
\begin{enumerate}
\item
If $f$ admits a right adjoint $g$, then $Ff$ admits the
right adjoint $Fg$.
\item
If $f$ is an equivalence in $\cA$, then $Ff$ is an equivalence
in $\cB$, and if $g$ is a quasi-inverse for $f$, then $Fg$ is
a quasi-inverse for $Ff$.
\item
If $F$ is fully faithful, then $f$ is an equivalence if and only
if $Ff$ is an equivalence in $\cB$.
\end{enumerate}
\end{lemma}
\begin{proof}(i): Indeed, let $\eps:f\circ g\Rightarrow\one_B$
and $\eta:\one_A\Rightarrow g\circ f$ be the unit and counit
of an adjunction for the pair $(f,g)$. We consider the $2$-cells
$$
\eps':=\delta^{F\ -1}_B\odot F(\eps)\odot\gamma^F_{g,f}:
Ff\circ Fg\Rightarrow\one_{FB}
\qquad
\eta':=\gamma^{F\ -1}_{f,g}\odot F(\eta)\odot\delta^F_A:
\one_{FA}\Rightarrow Fg\circ Ff
$$
and we compute :
$$
\begin{aligned}
(Fg*F\eps')\!\odot\!(F\eta'*Fg)=&\,(Fg*\delta^F_B)^{-1}\!\odot\!
(Fg*F\eps)\!\odot\!\gamma^{F\ -1}_{f\circ g,g}\!\odot\!\gamma^F_{g,g\circ f}
\!\odot\!(F\eta*Fg)\!\odot\!(\delta^F_A*Fg) \\
=&\,(Fg*\delta^F_B)^{-1}\odot\gamma^{F\ -1}_{\one_B,g}\odot
F(g*\eps)\odot F(\eta*g)\odot\gamma^F_{g,\one_B}\odot(\delta^F_A*Fg) \\
=&\,\one_{Fg}.
\end{aligned}
$$
Likewise we check that $(F\eps*Ff)\odot(Ff*F\eta)=\one_{Ff}$ : the
details are left to the reader.

(ii): It is clear from the proof of (i) that if $\eps$ and
$\eta$ are invertible $2$-cells, the same holds for $\eta'$
and $\eps'$, whence the contention.

(iii): Due to (ii), we may assume that $Ff$ is an equivalence
in $\cB$, and we check that $f$ is an equivalence in $\cA$.
To this aim, let $g:FB\to FA$ be a quasi-inverse for $Ff$,
and $\eta:\one_{FB}\to Ff\circ g$ and
$\eps:g\circ Ff\Rightarrow\one_{FA}$ two invertible $2$-cells.
Since $F_{AB}$ is an equivalence, we have a $1$-cell $h:B\to A$
and an invertible $2$-cell $\beta:Fh\Rightarrow g$, whence the
invertible $2$-cell
$$
\begin{aligned}
\eps':=\,&
\delta^F_A\odot\eps\odot(\beta*Ff)\odot\gamma^{F\ -1}_{f,h}:
F(h\circ f)\Rightarrow F\one_A \\
\eta':=\,&
\eta':=\,1
\gamma^F_{h,f}\odot(Ff*\beta^{-1})\odot\eta\odot\delta^{F\ -1}_B:
F\one_B\Rightarrow F(f\circ h).
\end{aligned}
$$
Again, since $F_{AB}$ is an equivalence, we deduce
invertible $2$-cells $\eps'':h\circ f\Rightarrow\one_A$
and $\eta'':\one_B\Rightarrow f\circ h$ such that
$F\eps''=\eps'$ and $F\eta''=\eta'$, whence the
contention.
\end{proof}

\begin{theorem}\label{th_pseudo-nat-equiv}
Let $\cA,\cB$ be two $2$-categories, $F,G:\cA\to\cB$ two
pseudo-functors, and $\alpha:F\Rightarrow G$ a pseudo-natural
transformation. The following conditions are equivalent :
\begin{enumerate}
\alphaenu
\item
$\alpha$ is a pseudo-natural equivalence of pseudo-functors.
\item
The $1$-cell $\alpha_A:FA\to GA$ is an equivalence in $\cB$,
for every $A\in\Ob(\cA)$.
\end{enumerate}
\end{theorem}
\begin{proof} Let $\cA^F$ be the $2$-category associated with
$F$ as in proposition \ref{prop_strictification}, together
with the strict pseudo-functors $F^\flat:\cA^F\to\cB$ and
$\pi^F:\cA^F\to\cA$, the pseudo-functor $\sigma^F:\cA\to\cA^F$
and the pseudo-natural isomorphism
$\omega^F:F\circ\pi^F\isom F^\flat$. Set
$\alpha':=(\alpha*\pi^F)\odot\omega^{F\ -1}:
F^\flat\Rightarrow G\circ\pi^F$, and suppose that there exists
a pseudo-natural transformation
$\beta':G\circ\pi^F\Rightarrow F^\flat$ with invertible
modifications $\Xi_1:\one_{F^\flat}\leadsto\beta'\odot\alpha'$
and $\Xi_2:\one_{G\circ\pi^F}\leadsto\alpha'\odot\beta'$. Then
notice that $\alpha'*\sigma^F=\alpha$, and set
$\beta:=\beta'*\pi^F$; we deduce invertible modifications
$\Xi_1\circ\pi^F:\one_F\leadsto\beta\odot\alpha$ and
$\Xi_2\circ\pi^F:\one_G\leadsto\alpha\odot\beta$. Thus,
we may replace $\cA$, $F$, $G$ and $\alpha$ by $\cA^F$,
$F^\flat$, $G\circ\pi^F$ and $\alpha'$, and assume from start
that $F$ is strict. Arguing likewise with $\cA^G$ and $G^\flat$,
we may then reduce to the case where both $F$ and $G$ are strict.

(b)$\Rightarrow$(a): Notice that the pseudo-functor
$\tilde\alpha:\cA\to 2\tdu\sMorph(\cB)$ associated with
$\alpha$ (see remark \ref{rem_pseudo-natural}(ii)) factors
through the inclusion pseudo-functor
$\Equiv(\cB)\to 2\tdu\sMorph(\cB)$ (notation of
\eqref{subsec_ups-and-equiv}); then define $\tilde\beta$
as the composition
$$
\cA\xrightarrow{\tilde\alpha}\Equiv(\cB)
\xrightarrow{\tilde\Upsilon}\Equiv(\cB)
$$
(notation of \eqref{subsec_tilde-upsilon}). Notice as
well that the coherence constraints of $\tilde\Upsilon$
are given by identities, and $\tilde\alpha$ is strict
(since $F$ and $G$ are strict), therefore the coherence
constraints of $\tilde\beta$ are given as well by identities,
so $\tilde\beta$ is strict. Furthermore, a simple inspection
shows that $\ss\circ\tilde\beta=G$ and $\st\circ\tilde\beta=F$,
where $\ss$ and $\st$ are the restrictions to $\Equiv(\cA)$ of
the source and target pseudo-functors of remark
\ref{rem_remark-cat-cats}(ii). We conclude that $\tilde\beta$
is the pseudo-functor associated as in remark
\ref{rem_pseudo-natural}(ii) with a pseudo-natural
transformation $\beta:G\Rightarrow F$. We can describe
the coherence constraint $\tau^\beta$ of $\beta$ as follows.
Let $\tau^\alpha$ be the coherence constraint of $\alpha$,
which assigns to every $1$-cell $f:A\to B$ of $\cA$ the
oriented square $\tau^\alpha_f$ as in remark
\ref{rem_pseudo-natural}(i). Recall that the pseudo-functor
$\sL$ of \eqref{subsec_ups-and-equiv} assigns to every
$1$-cell $g:X\to Y$ of $\cB$ a link
$\sL(g):=(g^\dagger,g,\eta^g,\eps^g):X\to Y$; then
$\tau^\alpha_f$ yields the $1$-cell
$\sL(\tau^\alpha_f):\sL(\alpha_A)\to\sL(\alpha_B)$ in the
$2$-category $\wLink(\cB)$. With this notation, we have
$\beta_A:=(\alpha_A)^\dagger$ for every $A\in\Ob(\cA)$, and
$\tau^\beta$ assigns to $f$ the oriented square
$$
\tau^\beta_f\qquad :\qquad
{\spreaddiagramcolumns{+80pt}\diagram
GA \ar[r]^-{\beta_A} \ar[d]_{Gf}
\drtwocell\omit{_\ \ \ \ \ \ \ \ \ \ \ \ \ \Upsilon(\sL(\tau^\alpha_f))^{-1}}
& FA \ar[d]^{Ff} \\
GB \ar[r]_-{\beta_B} & FB.
\enddiagram}
$$
To conclude, it remains to check that the systems of
$2$-cells :
$$
\lambda_A:=(\eps^{\alpha_A})^{-1}
\qquad\text{and}\qquad
\mu_A:=(\eta^{\alpha_A})^{-1}
\qquad
\text{for every $A\in\Ob(\cA)$}
$$
amount to two invertible modifications
$$
\lambda_\bullet:\one_F\leadsto\beta\odot\alpha
\qquad
\mu_\bullet:\alpha\odot\beta\leadsto\one_G.
$$
However, if $\tau^{\alpha\odot\beta}$ denotes the coherence
constraint of $\alpha\odot\beta$, the compatibility condition
for $\mu_\bullet$ comes down to the identity :
\set\begin{equation}\label{eq_comp-cond}
(\mu_B*Gf)\odot\tau^{\alpha\odot\beta}_f=Gf*\mu_A
\qquad
\text{for every $1$-cell $f:A\to B$ of $\cA$}.
\end{equation}
Let us then consider the diagram :
$$
\xymatrix@C+40pt{
GA \ar[r]^-{\beta_A} \ddouble \drtwocell\omit{_\ \ \ \ \one_{\beta_A}} &
FA \ddouble \rdouble \drtwocell\omit{_\ \ \ \ \one_{Ff}} &
FA \rdouble \ar[d]^{Ff} \drtwocell\omit{_\ \ \ \ \one_{Ff}} &
FA \ar[r]^-{\alpha_A} \ar[d]^{Ff} \drtwocell\omit{_\ \ \ \tau^\alpha_f} &
GA \ar[d]^{Gf} \\
GA \ar[r]|-{\beta_A} \ddouble \drtwocell\omit{_\ \ \ \mu_A} &
FA \ar[r]|-{Ff} \ar[d]^{\alpha_A}
\drtwocell\omit{_\ \ \ \ \ \ \ \ \ \sL(\tau^\alpha_f)^{-1}} &
FB \rdouble \ar[d]^{\alpha_B} \drtwocell\omit{_\ \ \ \lambda_B} &
FB \ar[r]|-{\alpha_B} \ddouble \drtwocell\omit{_\ \ \ \ \one_{\alpha_B}} &
GB \ddouble \\
GA \rdouble \ar[d]_{Gf} \drtwocell\omit{_\ \ \ \ \one_{Gf}} &
GA \ar[r]|-{Gf} \ar[d]^{Gf} \drtwocell\omit{_\ \ \ \ \one_{Gf}} &
GB \ar[r]|-{\beta_B} \ddouble \drtwocell\omit{_\ \ \ \mu_B} &
FB \ar[r]|-{\alpha_B} \ar[d]^{\alpha_B}
\drtwocell\omit{_\ \ \ \ \one_{\alpha_B}} &
GB \ddouble \\
GB \rdouble & GB \rdouble & GB \rdouble & GB \rdouble & GB.
}$$
We see that the composition of the squares of the top row
equals $\tau^\alpha_f*\beta_A$, and the composition of the
squares of the middle row equals $\alpha_B*\tau^\beta_f$.
Hence the composition of the squares of the top and middle
rows equals $\tau^{\alpha\odot\beta}_f$, and by further composing
with the squares of the bottom row we get the left hand-side
of \eqref{eq_comp-cond}. On the other hand, by proposition
\ref{prop_square-algebra} we can also compose the squares
column by column : then notice that the composition of the
squares in the third column from the left equals the identity
of $\alpha_B\circ Ff$. Thus, we may disregard this column,
and then notice that the composition of the squares of
the second column equals the inverse of the composition
of the squares of the fourth column. So we may disregard
all columns except the first from the left; but the composition
of the squares of the latter equals the right-hand side of
\eqref{eq_comp-cond}, as required. The verification for
$\lambda_\bullet$ is similar, and shall be left to the reader.

(a)$\Rightarrow$(b) is immediate from the definitions.
\end{proof}

\begin{corollary}\label{cor_first-corollary}
Let $\cA,\cB$ be two $2$-categories, $F,F':\cA\to\cB$
two pseudo-functors, and $\alpha:F\Rightarrow F'$ a
pseudo-natural equivalence. We have :
\begin{enumerate}
\item
For every $2$-category $\cA'$ and every pseudo-functor
$G:\cA'\to\cA$, the composition
$\alpha*G:F\circ G\Rightarrow F'\circ G$ is a pseudo-natural
equivalence (notation of remark {\em\ref{rem_pseudo-natural}(iv)}).
\item
For every $2$-category $\cB'$ and every pseudo-functor
$H:\cB\to\cB'$, the composition
$H*\alpha:H\circ F\Rightarrow H\circ F'$ is a pseudo-natural
equivalence.
\item
If $F':\cA\to\cB$ is any other pseudo-functor, and
$\alpha':F'\Rightarrow F''$ any other pseudo-natural
equivalence, then the composition
$\alpha'\odot\alpha:F\Rightarrow F''$ is a pseudo-natural
equivalence as well (notation of remark
{\em\ref{rem_pseudo-natural}(iii)}).
\end{enumerate}
\end{corollary}
\begin{proof} Assertions (i) and (ii) follow
immediately from theorem \ref{th_pseudo-nat-equiv}
and lemma \ref{lem_equiv-is-preserved}(ii). Assertion
(iii) follows from theorem \ref{th_pseudo-nat-equiv}
and lemma \ref{rem_compose-equiv}(ii) : details left
to the reader.
\end{proof}

\begin{definition}\label{def_2-equivalence}
Let $\cA$, $\cB$ be two $2$-categories, $F:\cA\to\cB$,
$G:\cB\to\cA$ two pseudo-functors. We say that $G$ is a
{\em right $2$-adjoint to $F$} (resp. a {\em strong right
$2$-adjoint to $F$}) if there exists a pseudo-natural
equivalence (resp. a pseudo-natural isomorphism) of
pseudo-functors
$$
\theta:H_\cB(F,\one_\cB)\Rightarrow H_\cA(\one_\cA,G)
$$
(definition \ref{def_using-mods}(i)). In this case,
we say that $(F,G)$ is a {\em $2$-adjoint pair}
(resp. a {\em strong $2$-adjoint pair}) of pseudo-functors,
and $\theta$ is a {\em $2$-adjunction} (resp. a {\em strong
$2$-adjunction}) for $(F,G)$.
\end{definition}

\begin{remark}\label{rem_2-adjoints}
According to theorem \ref{th_pseudo-nat-equiv}
and remark \ref{rem_pseudo-natural}(i), the datum of a
$2$-adjunction as in definition \ref{def_2-equivalence}
is equivalent to that of :
\begin{itemize}
\item
a system of equivalences of categories
$$
\theta_{AB}:\cB(FA,B)\to\cA(A,GB)
\qquad
\text{for every $A\in\Ob(\cA)$ and $B\in\Ob(\cB)$}
$$
\item
and for every pair of $1$-cells $f:A'\to A$ in $\cA$, $g:B\to B'$
in $\cB$, a natural isomorphism of functors
$$
\xymatrix{
\cB(FA,B) \ar[rrr]^-{\theta_{AB}} \ar[d]_{\cB(Ff,g)} &
\drtwocell\omit{_\ \ \ \ \ \ \ \tau^\theta_{(f,g)}} & &
\cA(A,GB) \ar[d]^{\cA(f,Gg)} \\
\cB(FA',B') \ar[rrr]_-{\theta_{A'B'}} & & & \cA(A',GB')
}$$
\item
such that, for every $2$-cells $\beta:f\Rightarrow f'$ in
$\cA$ and $\lambda:g\Rightarrow g'$ in $\cB$, we have the
identity
\set\begin{equation}\label{eq_deviancies}
\tau^\theta_{(f,g)}\boxminus H_\cB(F\beta,\lambda)=
H_\cA(\beta,G\lambda)\boxminus\tau^\theta_{(f',g')}
\end{equation}
\item
and for every composable pairs of $1$-cells
$A''\xrightarrow{f'}A'\xrightarrow{f}A$ in $\cA$, and
$B\xrightarrow{g}B'\xrightarrow{g'}B''$ in $\cB$, we
have the identities
$$
\begin{aligned}
(\tau^\theta_{(f',g')}\boxvert\tau^\theta_{(f,g)})\boxminus
H_\cB(\gamma^F_{f',f},\one_{g'\circ g})=\,&
H_\cA(\one_{f\circ f'},\gamma^G_{g,g'})\boxminus
\tau^\theta_{(f\circ f',g'\circ g)} \\
H_\cA(i_A,\delta^G_B)\boxminus\tau^\theta_{(\one_A,\one_B)}=\,&
\one_{\theta_{AB}}\boxminus H_\cB(\delta^F_A,i_B)
\end{aligned}
$$
where $(\delta^F,\gamma^F)$ (resp. $(\delta^G,\gamma^G)$)
denote by coherence constraint for $F$ (resp. for $G$).
\end{itemize}
\end{remark}

We wish to attach to every $2$-adjunction suitable units and
counits, as for usual adjoint pairs of functors; we shall see
that the triangular identities will have to be replaced by
certain invertible modifications. We begin with a few auxiliary
lemmata :

\begin{lemma}\label{lem_shouldbe-towards-2-yoneda}
Let $\cA$,$ \cB$ be any two $2$-categories, $F,F':\cA\to\cB$
two strict pseudo-functors, $\beta:F\Rightarrow F'$ a
pseudo-natural transformation. We have :
\begin{enumerate}
\item
$F$ induces a strict pseudo-natural transformation
(notation of example {\em\ref{ex_first-from-Cats-to-PsFuns}(ii)}) :
$$
H_F:H_\cA\Rightarrow H_\cB(F,F)
\qquad
(A,A')\mapsto(F_{AA'}:\cA(A,A')\to\cB(FA,FA')).
$$
\item
The coherence constraint  $\tau^\beta$ of $\beta$ induces
an invertible modification
$$
H_\beta:H_\cB(\beta,\one_{F'})\odot H_{F'}\leadsto
H_\cB(\one_F,\beta)\odot H_F
$$
that assigns to every object $(A_1,A_2)$ of $\cA^o\times\cA$
the natural transformation
$$
\tau^\beta_{A_1A_2}:H_\cB(\beta_{A_1},\one_{F'A_2})\circ F'_{A_1A_2}
\Rightarrow H_\cB(\one_{FA_1},\beta_{A_2})\circ F_{A_1A_2}.
$$
\item
If $F$ is fully faithful, $H_F$ is a pseudo-natural equivalence.
\end{enumerate}
\end{lemma}
\begin{proof} Assertion (i) is straightforward, and (iii)
follows immediately from theorem \ref{th_pseudo-nat-equiv}.

(ii): For every $1$-cell $(f_1,f_2):(A_1,A_2)\to(A'_1,A'_2)$
of $\cA^o\times\cA$ we consider the diagrams
$$
\xymatrix{ \cA(A_1,A_2) \ar[rrr]^-{H_{F'}} \ar[d]_{\cA(f_1,f_2)}
\drrtwocell\omit{_\ \ \one}
& & & \cB(F'A_1,F'A_2) \ar[rrr]^-{\cB(\beta_{A_1},\one_{F'A_2})}
\ar[d]|{\cB(F'f_1,F'f_2)} 
\drrtwocell\omit{_\qquad\qquad\quad\ \cB((\tau^\beta_{f_1})^{-1},\one_{F'f_2})}
& & & \cB(FA_1,F'A_2) \ar[d]^{\cB(Ff_1,F'f_2)} \\
\cA(A'_1,A'_2) \ar[rrr]_-{H_{F'}} \ddouble & &
\drtwocell\omit{_\qquad\tau^\beta_{A'_1A'_2}} &
\cB(F'A'_1,F'A'_2) \ar[rrr]_-{\cB(\beta_{A'_1},\one_{F'A'_2})}
& & & \cB(FA'_1,F'A'_2) \ddouble \\
\cA(A'_1,A'_2) \ar[rrr]_-{H_F} & & &
\cB(FA'_1,FA'_2) \ar[rrr]_-{\cB(\one_{FA'_1},\beta_{A_2})}
& & & \cB(FA'_1,F'A'_2)
}$$
$$
\xymatrix{ \cA(A_1,A_2) \ar[rrr]^-{H_{F'}}  \ddouble  & &
\drtwocell\omit{_\qquad\tau^\beta_{A_1A_2}} &
\cB(F'A_1,F'A_2) \ar[rrr]^-{\cB(\beta_{A_1},\one_{F'A_2})}
& & & \cB(FA_1,F'A_2) \ddouble \\
\cA(A_1,A_2) \ar[rrr]^-{H_F} \ar[d]_{\cA(f_1,f_2)}
\drrtwocell\omit{_\ \one}
& & & \cB(FA_1,FA_2) \ar[rrr]^-{\cB(\one_{FA_1},\beta_{A_2})}
\ar[d]|{\cB(Ff_1,Ff_2)} 
\drrtwocell\omit{_\qquad\qquad\cB(\one_{Ff_1},\tau^\beta_{f_2})}
& & & \cB(FA_1,F'A_2) \ar[d]^{\cB(Ff_1,F'f_2)} \\
\cA(A'_1,A'_2) \ar[rrr]_-{H_F} & & &
\cB(FA'_1,FA'_2) \ar[rrr]_-{\cB(\one_{FA'_1},\beta_{A'_2})}
& & & \cB(FA'_1,F'A'_2)
}$$
and notice that the coherence constraint of
$H_\cB(\beta,\one_{F'})\odot H_{F'}$ is given by the compositions
of oriented squares on the top row of the first diagram (see
example \ref{ex_from-Cats-to-PsFuns}). Likewise, the coherence
constraint of $H_\cB(\one_F,\beta)\odot H_F$ is given by the
composition of the oriented squares on the bottom row of the
second diagram. Therefore, we come down to checking that the
composition of the oriented squares of the first diagram equals
the composition of the oriented squares of the second diagram.
The latter translates as the identity :
$$
\tau^\beta_{f_2\circ t\circ f_1}\odot
(F'(f_2\circ t)*\tau^{\beta\ -1}_{f_1})=
(\tau^\beta_{f_2}*F(t\circ f_1))\odot(F'f_2*\tau^\beta_t*Ff_1)
$$
for every $1$-cell $t:A_1\to A_2$ in $\cA$. However, the
coherence axiom for $\beta$ yields the identities :
$$
\begin{aligned}
(\tau^\beta_t*Ff_1)\odot(F't*\tau^\beta_{f_1})&\,=
\tau^\beta_{t\circ f_1} \\
(\tau_{f_2}^\beta*F(t\circ f_1))\odot(F'f_2*\tau^\beta_{t\circ f_1})
&\,=\tau^\beta_{f_2\circ t\circ f_1}.
\end{aligned}
$$
The sought identity is an immediate consequence.
\end{proof}

\sset\subsubsection{}\label{subsec_converse-2-yoneda}
Conversely, consider now two strict pseudo-functors $F,G:\cA\to\cB$,
and a pseudo-natural transformation
$\lambda:H_\cA\Rightarrow H_\cB(F,G)$ with coherence constraint
$\tau_\bullet$. We set
$$
\lambda^\vee_A:=\lambda_{(A,A)}(\one_A):FA\to GA
\qquad
\text{for every $A\in\Ob(\cA)$}
$$
and for every $1$-cell $f:A\to A'$ in $\cA$ we let
$$
\tau^\vee_f:Gf\circ\lambda^\vee_A\Rightarrow\lambda^\vee_{A'}\circ Ff
$$
be the unique $2$-cell of $\cB$ that fits into the commutative
diagram
$$
\xymatrix{
Gf\circ\lambda^\vee_A \ddouble
\ar@{=>}[rrrr]^-{\tau^\vee_f} & & & &
\lambda^\vee_{A'}\circ Ff \ddouble \\
Gf\circ\lambda^\vee_A\circ F\one_A
\ar@{=>}[rr]^-{\tau_{(\one_A,f),\one_A}} & & \lambda_{(A,A')}(f) & &
G\one_{A'}\circ\lambda^\vee_{A'}\circ Ff.
\ar@{=>}[ll]_-{\tau_{(f,\one_{A'}),\one_{A'}}}
}$$

\begin{lemma}\label{lem_kernel-2-yoneda}
With the notation of \eqref{subsec_converse-2-yoneda}, the
following holds :
\begin{enumerate}
\item
The system of $1$-cells $(\lambda^\vee_A~|~A\in\Ob(\cA))$ defines
a pseudo-natural transformation
$$
\lambda^\vee:F\Rightarrow G
$$
with coherence constraint given by the system of\/ $2$-cells
$\tau^\vee_\bullet$.
\item
$H_F^\vee=\one_F$.
\item
Let $\mu:H_\cA\Rightarrow H_\cB(F,G)$ be another pseudo-natural
transformations, and $\Theta:\lambda\leadsto\mu$ any modification.
Then we have a modification $\Theta^\vee:\lambda^\vee\leadsto\mu^\vee$
given by the rule :
$$
A\mapsto
(\Theta^\vee_A:=\Theta_{(A,A),\one_A}:\lambda^\vee_A\Rightarrow\mu^\vee_A)
\qquad
\text{for every $A\in\Ob(\cA)$}.
$$
\end{enumerate}
\end{lemma}
\begin{proof}(i): Let $f,f':A\to A'$ be two $1$-cells in $\cA$;
the naturality of $\tau^\vee_\bullet$ amounts to the commutativity
of the diagram
$$
\xymatrix@C+10pt{
Gf\circ\lambda^\vee_A
\ar@{=>}[rr]^-{\tau_{(\one_A,f),\one_A}}
\ar@{=>}[d]_{G\alpha*\lambda^\vee_A} & &
\lambda_{(A,A')}(f) \ar@{=>}[d]_{\lambda_{(A,A')}(\alpha)} & &
\lambda^\vee_{A'}\circ Ff
\ar@{=>}[ll]_-{\tau_{(f,\one_{A'}),\one_{A'}}}
\ar@{=>}[d]^{\lambda^\vee_{A'}*F\alpha} \\
Gf'\circ\lambda^\vee_A
\ar@{=>}[rr]^-{\tau_{(\one_A,f'),\one_A}} & &
\lambda_{(A,A')}(f') & &
\lambda^\vee_{A'}\!\circ\!Ff'
\ar@{=>}[ll]_-{\tau_{(f',\one_{A'}),\one_{A'}}}
}$$
for every $2$-cell $\alpha:f\Rightarrow f'$ in $\cA$. However,
the commutativity of the two square subdiagrams follows by
applying the naturality of $\tau$ to the $2$-cells $(i_A,\alpha)$
and $(\alpha,i_{A'})$ of $\cA^o\times\cA$.

Let us check the coherence axioms for $\tau^\vee_\bullet$.
First, since $F$ and $G$ are strict, we need to show :
$$
\tau^\vee_{\one_A}=\one_{\lambda^\vee_A}
\qquad
\text{for every $A\in\Ob(\cA)$}
$$
which follows by a simple inspection.
Next, let $A\xrightarrow{\ f\ }A'\xrightarrow{\ g\ }A''$ be
a composable pair of $1$-cell of $\cA$; we need to show that
$$
(\tau^\vee_g*Ff)\odot(Gg*\tau^\vee_f)=\tau^\vee_{g\circ f}.
$$
Now, the coherence axiom for $\tau_\bullet$ yields the identity
$$
\tau_{(\one_A,g),f}\odot(Gg*\tau_{(\one_A,f),\one_A})=\tau_{(\one_A,g\circ f),\one_A}.
$$
Hence we are reduced to showing that
$$
(\tau^\vee_g*Ff)\odot(Gg*\tau_{(f,\one_{A'}),\one_{A'}}^{-1})=
\tau_{(g\circ f,\one_{A''}),\one_{A''}}^{-1}\odot\tau_{(\one_A,g),f}.
$$
Next, again by the coherence axiom for $\tau_\bullet$ we see that
$$
\tau_{(g\circ f,\one_{A''}),\one_{A''}}=
\tau_{(f,\one_{A''}),g}\odot(\tau_{(g,\one_{A''}),\one_{A''}}*Ff).
$$
Hence we are further reduced to checking the identity :
\set\begin{equation}\label{eq_anne-hathaway}
\tau_{(f,\one_{A''}),g}\odot(\tau_{(\one_{A'},g),\one_{A'}}*Ff)
=\tau_{(\one_A,g),f}\odot(Gg*\tau_{(f,\one_{A'}),\one_{A'}}).
\end{equation}
However, applying twice the coherence axiom for $\tau_\bullet$
to the identities in $\cA^o\times\cA$ :
$$
(f,\one_{A''})\circ(\one_{A'},g)=(f,g)=(\one_A,g)\circ(f,\one_{A'})
$$
we see that both sides of \eqref{eq_anne-hathaway} equal
$\tau_{(f,g),\one_{A'}}$.

(ii) follows by a direct inspection of the definitions.

(iii): The assertion follows from the commutativity of the diagram :
$$
\xymatrix{ Gf\circ\lambda^\vee_A
\ar@{=>}[rr]^-{\tau^\lambda_{(\one_A,f),\one_A}}
\ar@{=>}[d]_{Gf*\Theta^\vee_A} & &
\lambda_{(A,A')}(f) \ar@{=>}[d]_{\Theta_{(A,A'),f}} & &
\ar@{=>}[d]^{\Theta^\vee_{A'}*Ff}
\lambda^\vee_{A'}\circ Ff.
\ar@{=>}[ll]_-{\tau^\lambda_{(f,\one_{A'}),\one_{A'}}} \\
Gf\circ\mu^\vee_A
\ar@{=>}[rr]^-{\tau^\mu_{(\one_A,f),\one_A}} & & \mu_{(A,A')}(f) & &
\mu^\vee_{A'}\circ Ff
\ar@{=>}[ll]_-{\tau^\mu_{(f,\one_{A'}),\one_{A'}}}
}$$
which in turn follows directly from the definition of $\Theta$.
\end{proof}

\sset\subsubsection{}\label{subsec_2-units-for-2-adj}
Let now $F:\cA\to\cB$ be a pseudo-functor, and $G:\cB\to\cA$
a right $2$-adjoint for $F$, and $\theta$ a $2$-adjunction
for the pair $(F,G)$ as in definition \ref{def_2-equivalence}.
We suppose morever that $F$ and $G$ are strict; by lemma
\ref{lem_long-time}(i) we may find a pseudo-natural equivalence
$$
\psi:H_\cA(\one_\cA,G)\Rightarrow H_\cB(F,\one_\cB)
$$
as well as invertible modifications
$$
\Xi:\one_{H_\cB(F,\one_\cB)}\leadsto\psi\odot\theta
\qquad
\Theta:\theta\odot\psi\leadsto\one_{H_\cA(\one_\cA,G)}
$$
such that
$$
(\psi*\Theta)\odot(\Xi*\psi)=\one_\psi
\qquad\text{and}\qquad
(\Theta*\theta)\odot(\theta*\Xi)=\one_\theta.
$$
To ease notation, we shall drop the subscripts when referring
to $\theta_{AB}(h)$ or $\theta_{AB}(\mu)$ for any $1$-cell
$h:FA\to B$ or $2$-cell $\mu:h\Rightarrow h'$ in $\cB$, and
write simply $\theta(h)$ and $\theta(\mu)$. Likewise, we shall
write simply $\psi(k)$ and $\psi(\nu)$ for any $1$-cell
$k:A\to GB$ and any $2$-cell $\nu:k\Rightarrow k'$ in $\cA$.
Now, for every $A\in\Ob(\cA)$ and $B\in\Ob(\cB)$ let us set
$$
\eta_A:=\theta(\one_{FA}):A\to GFA
\qquad
\eps_B:=\psi(\one_{GB}):FGB\to B.
$$
For every $1$-cell $f:A\to A'$ in $\cA$, and every $1$-cell
$g:B\to B'$ in $\cB$ we let
$$
\tau^\eta_f:GFf\circ\eta_A\Rightarrow\eta_{A'}\circ f
\qquad\text{and}\qquad
\tau^\eps_g:g\circ\eps_B\Rightarrow\eps_{B'}\circ FGg
$$
be the unique $2$-cells in $\cA$ and respectively in $\cB$
that make commute the diagrams :
$$
\xymatrix{
GFf\circ\eta_A \ar@{=>}[rd]_{\tau^\theta_{(\one_A,Ff),\one_{FA}}}
\ar@{=>}[rr]^-{\tau^\eta_f} & &
\eta_{A'}\circ f \ar@{=>}[ld]^{\ \ \tau^\theta_{(f,\one_{FA'}),\one_{FA'}}} &
g\circ\eps_B \ar@{=>}[rr]^-{\tau^\eps_g}
\ar@{=>}[rd]_{\tau^\psi_{(\one_{GB},g),\one_{GB}}}
& & \eps_{B'}\circ FGg \ar@{=>}[ld]^{\ \ \tau^\psi_{(Gg,\one_{B'}),\one_{GB'}}} \\
& \theta(Ff) & & & \psi(Gg)
}$$
where, as usual, $\tau^\psi$ denotes the coherence constraint
of $\psi$.

\begin{lemma}\label{lem_2-triangular-ident}
With the notation of \eqref{subsec_2-units-for-2-adj}, we have :
\begin{enumerate}
\item
The system $(\eta_A~|~A\in\Ob(\cA))$
defines a pseudo-natural transformation
$$
\eta:\one_\cA\Rightarrow GF
$$
whose coherence constraint is given by the system of\/ $2$-cells
$\tau^\eta_\bullet$.
\item
The system $(\eps_B~|~B\in\Ob(\cB))$ defines a pseudo-natural
transformation
$$
\eps:FG\Rightarrow\one_\cB
$$
whose coherence constraint is given by the system of\/ $2$-cells
$\tau^\eps_\bullet$.
\end{enumerate}
\end{lemma}
\begin{proof}(i): Define the pseudo-natural transformation
$H_F:H_\cA\Rightarrow H_\cB(F,F)$ as in lemma
\ref{lem_shouldbe-towards-2-yoneda}(i), and set
$$
\lambda:=(\theta*(\one_{\cA^o}\times F))\odot H_F:
H_\cA\Rightarrow H_\cA(\one_\cA,G\circ F).
$$
Lemma \ref{lem_kernel-2-yoneda}(i) yields a pseudo-natural
transformation
$$
\lambda^\vee:\one_\cA\Rightarrow G\circ F
$$
and it is easily seen that $\lambda^\vee_A=\eta_A$ for every
$A\in\Ob(\cA)$. It remains to check that $\tau^\eta$ agrees
with the coherence constraint $\tau^{\lambda^\vee}$ of
$\lambda^\vee$. However, denote by $\tau^\lambda$, $\tau^{H_F}$
and $\tau^{\theta*(\one_\cA\times F)}$ the coherence constraints of
respectively $\lambda$, $H_F$ and $\theta*(\one_\cA\times F)$;
by inspecting the definitions we find that
$$
\tau^\lambda_{(\one_A,f),\one_A}=
(\theta*(\one_\cA\times F))_{A,A'}(\tau^{H_F}_{(\one_A,f),\one_A})
\odot\tau^{\theta*(\one_\cA\times F)}_{(\one_A,f),\one_{FA}}=
\tau^\theta_{(\one_A,Ff),\one_{FA}}
$$
and likewise : $\tau^\lambda_{(f,\one_{A'}),\one_{A'}}=
\tau^\theta_{(f,\one_{FA'}),\one_{FA'}}$. The assertion follows.

(ii): A similar calculation shows that
$\eps=((\psi*(G\times\one_\cB))\odot H_G)^\vee$.
\end{proof}

We now turn to the case of an arbitrary $2$-adjoint pair
of pseudo-functors :

\begin{theorem}\label{th_2-adjunction}
Let $F:\cA\to\cB$ and $G:\cB\to\cA$ be two pseudo-functors.
We have :
\begin{enumerate}
\item
If $(F,G)$ is a $2$-adjoint pair, there exist pseudo-natural
transformations
$$
\eta:\one_\cA\Rightarrow GF
\qquad
\eps:FG\Rightarrow\one_\cB
$$
related by a pair of invertible modifications
$$
\Sigma:(G*\eps)\odot(\eta*G)\leadsto\one_G
\qquad
\Sigma':(\eps*F)\odot(F*\eta)\leadsto\one_F.
$$
We call $\eta$ a {\em unit} and $\eps$ a {\em counit} for
the $2$-adjunction $\theta$, and $\Sigma$ and $\Sigma'$ the
{\em triangular modifications} associated with $(\eta,\eps)$.
\item
Conversely, the existence of pseudo-natural transformations
$\eps$, $\eta$ and invertible modifications $\Sigma$, $\Sigma'$
as in {\em (i)}, implies that $G$ is right $2$-adjoint to $F$.
\end{enumerate}
\end{theorem}
\begin{proof} After replacing $\sU$ by a larger universe, we
may assume that $\cA$ and $\cB$ are small; in this case, let
us remark that the existence of pseudo-natural transformations
$\eta$ and $\eps$ related by triangular modifications as in
(i) means precisely that $(F,G)$ is an adjoint pair of $1$-cells
of the $2$-category $\overline{2\tdu\bCat}$ of remark
\ref{rem_reduced-2-cats}.
Let $F^u$ and $G^u$ be the unital pseudo-functors associated
with $F$ and $G$, with pseudo-natural isomorphisms
$\alpha^F:F\Rightarrow F^u$ and $\alpha^G:G\Rightarrow G^u$,
as in proposition \ref{prop_towards-2-yoneda}. If $\theta$
is a $2$-adjunction for $(F,G)$, we get a unique pseudo-natural
equivalence $\theta^*$ fitting into the commutative diagram :
$$
\xymatrix{ H_\cB(F,\one_\cB)
\ar@{=>}[d]_{H_\cB(\alpha^F,i_\cB)} \ar@{=>}[r]^-\theta &
H_\cA(\one_\cA,G) \ar@{=>}[d]^{H_\cA(i_\cA,\alpha^G)} \\
H_\cB(F^u,\one_\cB) \ar@{=>}[r]^-{\theta^*} &
H_\cA(\one_\cA,G^u)
}$$
(notation of example \ref{ex_first-from-Cats-to-PsFuns}(ii);
here $i_\cA$ is the identity transformation of the
pseudo-functor $\one_\cA$, and likewise for $i_\cB$).
Then $\theta^*$ is a $2$-adjunction for the pair
$(F^u,G^u)$. By the same token, from a given $2$-adjunction
for the pair $(F^u,G^u)$ we easily deduce a $2$-adjunction
for $(F,G)$. Also, $(F,G)$ is an adjoint pair of $1$-cells
in $\overline{2\tdu\bCat}$ if and only if the same holds
for the pair $(F^u,G^u)$ (remark \ref{rem_compose-equiv}(iv)).
Thus, we may replace $F$ and $G$ by $F^u$ and $G^u$, and assume
from start that $F$ and $G$ are unital pseudo-functors.
Next we consider the $2$-categories $\cA^{F,G}$ and $\cB^{G,F}$
of proposition \ref{prop_strictify-pairs}, as well as the
strict pseudo-functors $F^\flat,G^\flat,\pi^{F,G}$ and $\pi^{G,F}$,
the pseudo-functors $\sigma^{F,G}$ and $\sigma^{G,F}$, and the
isomorphisms of pseudo-functors
$$
\psi^{F,G}:\sigma^{F,G}\circ\pi^{F,G}\isom\one_{\cA^{F,G}}
\qquad
\psi^{G,F}:\sigma^{G,F}\circ\pi^{G,F}\isom\one_{\cB^{G,F}}.
$$

\begin{claim}\label{cl_these-are-equivs}
$\pi^{F,G}$ and $\pi^{G,F}$ are $2$-equivalences.
\end{claim}
\begin{pfclaim} Since $\pi^{F,G}$ is the identity on objects,
it suffices to check that it is fully faithful. The latter
is clear, since $\pi^{F,G}\circ\sigma^{F,G}=\one_\cA$ and
$\sigma^{F,G}\circ\pi^{F,G}$ is isomorphic to $\one_{\cA^{F,G}}$
(details left to the reader). The same argument applies to
$\pi^{G,F}$.
\end{pfclaim}

Now, suppose that $\theta:H_\cB(F,\one_\cB)\isom H_\cA(\one_\cA,G)$
is a $2$-adjunction for the pair $(F,G)$. By claim
\ref{cl_these-are-equivs} and lemma
\ref{lem_shouldbe-towards-2-yoneda}(iii), $\pi^{F,G}$ and
$\pi^{G,F}$ induce pseudo-natural equivalences
$$
H_{\pi^{F,G}}:H_{\cA^{F,G}}\isom H_\cA(\pi^{F,G},\pi^{F,G})
\qquad\text{and}\qquad
H_{\pi^{G,F}}:H_{\cB^{G,F}}\isom H_\cB(\pi^{G,F},\pi^{G,F}).
$$
Then we may find a pseudo-natural equivalence
$\theta^\flat:H_{\cB^{G,F}}(F^\flat,\one_{\cB^{G,F}})\isom
H_{\cA^{F,G}}(\one_{\cA^{F,G}},G^\flat)$ fitting into
the essentially commutative diagram of pseudo-functors :
$$
\xymatrix@C+20pt{ H_{\cB^{G,F}}(F^\flat,\one_{\cB^{G,F}})
\ar[d]_{H_{\pi^{G,F}}*(F^{\flat o}\times\one_{\cB^{G,F}})}
\ar[rr]^-{\theta^\flat} & & H_{\cA^{F,G}}(\one_{\cA^{F,G}},G^\flat)
\ar[d]^{H_{\pi^{F,G}}*(\one_{\cA^{F,G}}\times G^\flat)} \\
H_\cB(\pi^{G,F}\circ F^\flat,\pi^{G,F})
\ar[d]_{H_\cB*(\omega^{F,G}\times\one_{\pi^{G,F}})} & &
H_\cA(\pi^{F,G},\pi^{F,G}\circ G^\flat) \\
H_\cB(F\circ\pi^{F,G},\pi^{G,F})
\ar[rr]^-{\theta*(\pi^{F,G}\times\pi^{G,F})} & &
\ar[u]_{H_\cA*(\one_{\pi^{F,G}}\times\omega^{G,F})}
H_\cA(\pi^{F,G},G\circ\pi^{G,F}).
}$$
Conversely, if $\theta^\flat$ is a given $2$-adjunction
for the pair $(F^\flat,G^\flat)$, then we deduce a pseudo-natural
equivalence $\theta$ fitting into the essentially commutative
diagram :
$$
\xymatrix@C+20pt{ H_\cB(F,\one_\cB)
\ar[rr]^-\theta & & H_\cA(\one_\cA,G) \\
H_{\cB^{G,F}}(\sigma^{G,F}\circ F,\sigma^{G,F})
\ar[u]^{H_{\pi^{G,F}}*((\sigma^{G,F}\circ F)^o\times\sigma^{G,F})}
\ddouble & & H_{\cA^{F,G}}(\sigma^{F,G},\sigma^{F,G}\circ G)
\ar[u]_{H_{\pi^{F,G}}*(\sigma^{F,G}\times(\sigma^{F,G}\circ G))} \\
H_{\cB^{G,F}}(F^\flat\circ\sigma^{F,G},\sigma^{G,F})
\ar[rr]^-{\theta^\flat*(\sigma^{F,G}\times\sigma^{G,F})} & &
\udouble H_{\cA^{F,G}}(\sigma^{F,G},G^\flat\circ\sigma^{G,F}).
}$$
Thus $(F,G)$ is a $2$-adjoint pair if and only if the same
holds for $(F^\flat,G^\flat)$.

Likewise, if $\eta:\one_\cA\Rightarrow GF$ and
$\eta:FG\Rightarrow\one_\cB$ are given pseudo-natural
transformations fulfilling the triangular identities
$(G*\eps)\odot(\eta*G)=\one_G$ and
$(\eps*F)\odot(F*\eta)=\one_F$ in $\overline{2\tdu\bCat}$,
there follow pseudo-natural transformations :
$$
\xymatrix@C+30pt@R-15pt{
\eta^\flat:\one_{\cA^{F,G}} \ar@{=>}[r]^-{\psi^{F,G\,-1}} &
\sigma^{F,G}\circ\pi^{F,G} \ar@{=>}[r]^-{\sigma^{F,G}*\eta*\pi^{F,G}}
& \sigma^{F,G}\circ GF\circ\pi^{F,G} \ddouble \\
G^\flat F^\flat & G^\flat F^\flat\circ\sigma^{F,G}\circ\pi^{F,G}
\rdouble \ar@{=>}[l]_-{G^\flat F^\flat*\psi^{F,G}} &
G^\flat\circ\sigma^{G,F}\circ F\circ\pi^{F,G} \\
\eps^\flat:F^\flat G^\flat \ar@{=>}[r]^-{F^\flat G^\flat*\psi^{G,F\,-1}}
& F^\flat G^\flat\circ\sigma^{G,F}\circ\pi^{G,F} \rdouble &
F^\flat\circ\sigma^{F,G}\circ G\circ\pi^{G,F} \ddouble \\
\one_{\cB^{F,G}} & \sigma^{G,F}\circ\pi^{G,F} \ar@{=>}[l]^-{\psi^{G,F}}
& \sigma^{G,F}\circ FG\circ\pi^{G,F}
\ar@{=>}[l]^-{\sigma^{G,F}*\eps*\pi^{G,F}}
}$$
Set $X:=G^\flat*\psi^{G,F}$ and $Y:=\psi^{F,G}*G^\flat$.
By proposition \ref{prop_strictify-pairs}(ii) we have :
$$
Z:=G^\flat F^\flat*((G^\flat*\psi^{G,F})^{-1}\odot(\psi^{F,G}*G^\flat))=
(G^\flat F^\flat\sigma^{F,G}*\omega^{G,F})^{-1}=
(\sigma^{F,G}GF*\omega^{G,F})^{-1}
$$
whence : $Z\odot(\sigma^{F,G}*\eta*\pi^{F,G}*G^\flat)=
\sigma^{F,G}*((\eta*G\pi^{G,F})\odot\omega^{G,F\,-1})$.
We deduce :
$$
\begin{aligned}
(G^\flat*\eps^\flat)\odot(\eta^\flat*G^\flat)
&\,=X\odot(G^\flat\sigma^{G,F}*\eps*\pi^{G,F})
\odot Z\odot(\sigma^{F,G}*\eta*\pi^{F,G}*G^\flat)\odot Y^{-1} \\
&\,=X\odot\!(G^\flat\sigma^{G,F}*\eps*\pi^{G,F})\odot
(\sigma^{F,G}*((\eta*G\pi^{G,F})\odot\omega^{G,F\,-1}))\odot Y^{-1} \\
&\,=X\odot(\sigma^{F,G}G*\eps*\pi^{G,F})\odot
(\sigma^{F,G}*((\eta*G\pi^{G,F})\odot\omega^{G,F\,-1}))\odot Y^{-1} \\
&\,=X\odot(\sigma^{F,G}*\omega^{G,F\,-1})\odot Y^{-1} \\
&\,=\one_{G^\flat}
\end{aligned}
$$
where the last identity follows again from proposition
\ref{prop_strictify-pairs}(ii). A similar computation that
we leave to the reader shows that
$(\eps^\flat*F^\flat)\odot(F^\flat*\eta^\flat)=\one_{F^\flat}$ in
$\overline{2\tdu\bCat}$. Conversely, if given pseudo-natural
transformations $\eta^\flat:\one_{\cA^{F,G}}\Rightarrow G^\flat F^\flat$
and $\eps^\flat:F^\flat G^\flat\Rightarrow\one_{\cB^{G,F}}$ satisfy
the triangular identities, then we get pseudo-natural
transformations
$$
\eta:=\pi^{F,G}*\eta^\flat*\sigma^{F,G}:\one_\cA\Rightarrow GF
\qquad
\eps:=\pi^{G,F}*\eps^\flat*\sigma^{G,F}:FG\Rightarrow\one_\cB
$$
that also fulfill the corresponding triangular identities
(details left to the reader). Thus, we may replace $F$ and
$G$ by $F^\flat$ and $G^\flat$, and assume from start that
$F$ and $G$ are strict.

In order to exhibit $\Sigma$, we mimick the discussion of
\eqref{subsec_adj-pair} : first, for every $B\in\Ob(\cB)$
we have an invertible $2$-cell
$$
{\diagram GB \ar[rrrr]^-{(\eta*G)_B}
\ar@/_1pc/[rrrrd]_{\theta\odot\psi(\one_{GB})} &
\drrtwocell\omit{\qquad\qquad\qquad\tau^\theta_{(\one_{GB},\eps_B),\one_{FGB}}}
& & &
GFGB \ar[d]^{(G*\eps)_B} \\
& & & & GB 
\enddiagram}
$$
where $\tau^\theta$ denotes the coherence constraint of
$\theta$, and $\psi$ is the pseudo-natural equivalence as
in \eqref{subsec_2-units-for-2-adj} used to define $\eps$.
Then, with the notation of \eqref{subsec_2-units-for-2-adj},
we set
$$
\Sigma_B:=\Theta_{GB,B}\odot\tau^\theta_{(\one_{FGB},\eps_B),\one_{FGB}}
\qquad
\text{for every $B\in\Ob(\cB)$}.
$$
We have to check that the system $(\Sigma_B~|~B\in\Ob(\cB))$
is a modification $(G*\eps)\odot(\eta*G)\leadsto\one_G$.
However, set (notation of remark
\ref{rem_construct-modifs}(iv,v) and lemma
\ref{lem_shouldbe-towards-2-yoneda}(i))
$$
\mu:=((\theta\odot\psi)*(G\times\one_\cB))\odot H_G
\qquad\text{and}\qquad
\Theta':=(\Theta\circ(G\times\one_\cB))*H_G:
\mu\leadsto H_G.
$$
Unwinding the definitions, we see that
$\Theta'_{B,B'}=\Theta_{GB,B'}*G_{BB'}$ for every
$(B,B')\in\Ob(\cB^o\times\cB)$; according to lemma
\ref{lem_kernel-2-yoneda}(ii,iii) we then obtain an
invertible modification
$$
\Theta'^\vee:\mu^\vee\leadsto\one_G
\qquad
B\mapsto\Theta_{GB,B}.
$$
We are thus reduced to showing that the system
$(\tau^\theta_{(\one_{GB},\eps_B),\one_{FGB}}~|~B\in\Ob(\cB))$
yields an invertible modification
$$
\Theta'':(G*\eps)\odot(\eta*G)\leadsto\mu^\vee
$$
since in this case we shall have
$\Sigma=\Theta'^\vee\odot\Theta''$. However, by unwinding
the definitions we see that the coherence constraint of
$\mu^\vee$ is the system of $2$-cells
$$
\tau^\mu_g:=(\tau^\theta_{(Gg,\one_{B'}),\eps_{B'}})^{-1}\odot
\theta(\tau^\eps_g)\odot\tau^\theta_{(\one_{GB},g),\eps_B}
\qquad
\text{for every $1$-cell $g:B\to B'$ of $\cB$}
$$
whereas that of $(G*\eps)\odot(\eta*G)$ is given by the
system of $2$-cells
$$
\tau^{(G*\eps)\odot(\eta*G)}_g:=(G(\eps_{B'})*\tau^\eta_{Gg})
\odot(G(\tau^\eps_g)*\eta_{GB})
\qquad
\text{for every $g:B\to B'$}
$$
and we need to show the identities :
$$
(\tau^\theta_{(\one_{GB'},\eps_{B'}),\one_{FGB'}}*Gg)\odot
\tau^{(G*\eps)\odot(\eta*G)}_g=X:=\tau^\mu_g\odot
(Gg*\tau^\theta_{(\one_{FGB},\eps_B),\one_{FGB}})
$$
for every $1$-cell $g:B\to B'$ of $\cB$. We compute :
$$
\begin{aligned}
X=&\,
(\tau^\theta_{(Gg,\one_{B'}),\eps_{B'}})^{-1}\odot\theta(\tau^\eps_g)
\odot\tau^\theta_{(\one_{GB},g\circ\eps_B),\one_{FGB}} \\
=&\, (\tau^\theta_{(Gg,\one_{B'}),\eps_{B'}})^{-1}\odot
\tau^\theta_{(\one_{GB},\eps_{B'}\circ FGg),\one_{FGB}}\odot
(G(\tau^\eps_g)*\eta_{GB}) 
\end{aligned}
$$
where the first equality follows from the coherence axioms for
$\tau^\theta$, and the second follows from \eqref{eq_deviancies}.
So, we are further reduced to checking the identities :
$$
\tau^\theta_{(Gg,\one_{B'}),\eps_{B'}}\odot
(\tau^\theta_{(\one_{GB'},\eps_{B'}),\one_{FGB'}}*Gg)\odot
(G(\eps_{B'})*\tau^\eta_{Gg})=
\tau^\theta_{(\one_{GB},\eps_{B'}\circ FGg),\one_{FGB}}.
$$
However, by applying again the coherence axiom for
$\tau^\theta$ we see that
$$
\begin{aligned}
\tau^\theta_{(\one_{GB},\eps_{B'}\circ FGg),\one_{FGB}}=&\,
\tau^\theta_{\one_{GB},\eps_{B'},FGg}\odot
(G(\eps_{B'})*\tau^\theta_{(\one_{GB},FGg),\one_{FGB}}) \\
\tau^\theta_{(Gg,\one_{B'}),\eps_{B'}}\odot
(\tau^\theta_{(\one_{GB'},\eps_{B'}),\one_{FGB'}}*Gg)=&\,
\tau^\theta_{(Gg,\eps_{B'}),\one_{FGB'}}
\end{aligned}
$$
so we are further reduced to showing :
$$
\tau^\theta_{(Gg,\eps_{B'}),\one_{FGB'}}=
\tau^\theta_{\one_{GB},\eps_{B'},FGg}\odot
(G(\eps_{B'})*\tau^\theta_{(Gg,\one_{FGB'}),\one_{FGB'}}).
$$
But the latter follows by yet another application of
the coherence axioms for $\tau^\theta$. This concludes
the construction of $\Sigma$. Concerning $\Sigma'$,
notice that for every $A\in\Ob(\cA)$ and $B\in\Ob(\cB)$
$$
\cA^o(G^oB^o,A^o)=\cA(A,GB)
\qquad\text{and}\qquad
\cB^o(B^o,F^oA^o)=\cB(FA,B)
$$
and $\psi$ can thus be regarded as a system of equivalences of
categories $\psi^o_{B^oA^o}:\cA^o(G^oB^o,A^o)\to\cB^o(B^o,F^oA^o)$.
The latter amounts then to a pseudo-natural equivalence
$$
\psi^o:\cA^o(G^o,\one_{\cA^o})\Rightarrow\cB^o(\one_{\cB^o},F^o)
$$
with coherence constraint given by the system of isomorphisms of
functors $\tau^{\psi^o}_{(g^o,f^o)}:=\tau^\psi_{(f,g)}$, for every
$2$-cell $f$ of $\cA$ and $g$ of $\cB$. Likewise, we may define
the $2$-adjunction $\theta^o$ for the pair $(G^o,F^o)$, and a
simple inspection shows that performing the foregoing
constructions on this new pair of pseudo-natural equivalences
amounts to swapping the roles of $\eps$ and $\eta$; then the
modification $\Sigma^o$ for the new pair of unit and counit
will give the sought $\Sigma'$ for the original pair.

(ii): Given such $\eta$ and $\eps$, we set
$$
\begin{aligned}
\theta:=&\,
H_\cA(\eta,\one_G)\odot(H_G*(F^o\times\one_\cB)):
H_\cB(F,\one_\cB)\Rightarrow H_\cA(\one_\cA,G) \\
\psi:=&\,
H_\cB(\one_F,\eps)\odot(H_F*(\one^o_\cA\times G)):
H_\cA(\one_\cA,G)\Rightarrow H_\cB(F,\one_\cB).
\end{aligned}
$$
Now, example \ref{ex_modifications}(ii) yields the
following identity in $\overline{2\tdu\bCat}$ :
$$
\Xi:(H_F*(\one^o_\cA\times G))\odot H_\cA(\eta,\one_G)
=(H_\cB(F,F)*(\eta^o\times\one_G))\odot(H_F*((GF)^o\times G))
$$
and on the other hand, it is easily seen that
$$
\begin{aligned}
(H_F\!*\!((GF)^o\times G))\odot(H_G\!*\!(F^o\times\one_\cB))
\!=\!H_{FG}\!*\!(F^o\times\one_\cB):&\,
H_\cB(F,\one_\cB)\!\Rightarrow\!H_\cB(FGF,FG) \\
(H_\cB(\one_F,\eps)\odot
(H_\cB(F,F)*(\eta^o\!\times\!\one_G))=
H_\cB(F*\eta,\eps):&\,
H_\cB(FGF,FG)\!\Rightarrow\!H_\cB(F,\one_\cB)
\end{aligned}
$$
whence the identity in $\overline{2\tdu\bCat}$ :
$$
\psi\odot\theta=H_\cB(F*\eta,\eps)\odot
(H_{FG}*(F^o\times\one_\cB))=H_\cB(F*\eta,i_\cB)\odot
((H_\cB(\one_{FG},\eps)\odot H_{FG})*(F^o\times\one_\cB)).
$$
But we have $H_\cB(\one_{FG},\eps)\odot H_{FG}=H_\cB(\eps,i_\cB)$
in $\overline{2\tdu\bCat}$, by lemma
\ref{lem_shouldbe-towards-2-yoneda}(ii), so that :
$$
\begin{aligned}
\psi\odot\theta&\,=H_\cB(F*\eta,i_\cB)\odot
(H_\cB(\eps,i_\cB)*(F^o\times\one_\cB)) \\
&\,=H_\cB(F*\eta,i_\cB)\odot H_\cB(\eps*F,i_\cB) \\
&\,=H_\cB((\eps*F)\odot(F*\eta),i_\cB) \\
&\,=H_\cB(\one_F,i_\cB) \\
&\,=\one_{H_\cB(F,\one_\cB)}.
\end{aligned}
$$
Arguing likewise, we obtain as well the identity
$\theta\odot\psi=\one_{H_\cA(\one_\cA,G)}$ in
$\overline{2\tdu\bCat}$ : the details shall be left to
the reader. Summing up, we have shown that $\theta$ is
an equivalence in the $2$-category
$\sPsFun(\cA^o\times\cB,\bCat)$ with quasi-inverse $\psi$,
and the proof of the theorem is concluded.
\end{proof}

\begin{remark}\label{rem_opposing-thumb}
(i)\ \
Let $(F:\cA\to\cB,G:\cB\to\cA)$ be a $2$-adjoint pair of
pseudo-functors, and denote by $\eta$ and $\eps$ the unit
and the counit of a $2$-adjunction for this pair, and
$\Sigma$,$\Sigma'$ the associated triangular modifications
as in theorem \ref{th_2-adjunction}. There follow pseudo-natural
transformations
$$
\eta^o:G^oF^o\Rightarrow\one_{\cA^o}
\qquad
\eps^o:\one_{\cB^o}\Rightarrow F^oG^o
$$
and invertible modifications :
$$
\Sigma^o:(\eta^o*G^o)\odot(G^o*\eps^o)\leadsto\one_{G^o}
\qquad
\Sigma'^o:(F^o*\eta^o)\odot(\eps^o*F^o)\leadsto\one_{F^o}.
$$
In light of theorem \ref{th_2-adjunction}, it follows that
$(G^o,F^o)$ is a $2$-adjoint pair of pseudo-functors, and
$\eps^o$ and $\eta^o$ are respectively a unit and a counit
of a $2$-adjunction for this pair.

(ii)\ \
Likewise, we see that $({}^oF,{}^oG)$ is a $2$-adjoint pair
of pseudo-functors, and ${}^o\eta$ and ${}^o\eps$ are the
unit and respectively the counit of a $2$-adjunction for this
pair.

(iii)\ \
Moreover, the right $2$-adjoint $G$ of $F$ is unique up to
pseudo-natural equivalence of pseudo-functors. Indeed, we
have already observed that the pair $(F,G)$ is $2$-adjoint
if and only if it is an adjoint pair of $1$-cells of the
$2$-category $\sU'\tdu\overline{2\tdu\bCat}$ (for a suitable
universe $\sU'$); on the other hand, the pseudo-natural
equivalences of pseudo-functors are precisely the invertible
$2$-cells of $\sU'\tdu\overline{2\tdu\bCat}$, hence the
assertion follows from remark \ref{rem_compose-equiv}(v).
\end{remark}

\begin{corollary}\label{cor_fully-faith-2-adjoint}
Let $F:\cA\to\cB$ be a pseudo-functor and $G:\cB\to\cA$ a right
$2$-adjoint for $F$. Let also $\eta:\one_\cA\Rightarrow GF$ and
$\eps:FG\Rightarrow\one_\cB$ be the unit and counit of a
$2$-adjunction $\theta$ for the pair $(F,G)$. Then the following
conditions are equivalent :
\begin{enumerate}
\alphaenu
\item
$F$ (resp. $G$) is fully faithful.
\item
$\eta$ (resp. $\eps$) is a pseudo-natural equivalence.
\item
There exists a pseudo-natural equivalence $\eta':\one_\cA\isom GF$
(resp. $\eps':FG\isom\one_\cB$).
\end{enumerate}
Moreover, if\/ {\em (c)} holds, there exists a pseudo-natural
transformation $\eps'$ (resp. $\eta'$) with invertible
modifications $(G*\eps')\odot(\eta'*G)\leadsto\one_G$
and $(\eps'*F)\odot(F*\eta')\leadsto\one_F$.
\end{corollary}
\begin{proof} By inspecting the construction of $\eta$ and
$\eps$, we are easily reduced to the case where $\cA$ and $\cB$
are small and $F$ and $G$ are strict. Also, by remark
\ref{rem_opposing-thumb}(i) it suffices to check the assertions
for $F,\eta$ and $\eta'$. Suppose now that (a) holds for $F$.
Then the composition
$$
\lambda_{A'A}:\cA(A',A)\xrightarrow{F_{A'A}}\cB(FA',FA)
\xrightarrow{\theta_{A',FA}}\cA(A',GFA)
\qquad
\text{for every $A,A'\in\Ob(\cA)$}
$$
is an equivalence of categories, and $\eta_A=\lambda_{AA}(\one_A)$
for every $A\in\Ob(\cA)$. Taking $A':=GFA$, we find a $1$-cell
$f:GFA\to A$ with an invertible $2$-cell
$\lambda_{GFA,A}(f)\isom\one_{GFA}$. Moreover, as explained in the
proof of lemma \ref{lem_2-triangular-ident}(i), the rule :
$(A',A)\mapsto\lambda_{A'A}$ is pseudo-natural in both $A$ and $A'$;
thus, for every $1$-cell $h:A'\to A$ of $\cA$, the coherence
constraint $\tau^\lambda$ of $\lambda$ yields an isomorphism
of functors :
$$
\tau^\lambda_{(h,\one_A)}:\cA(h,\one_{GFA})\circ\lambda_{AA}\isom
\lambda_{GFA,A}\circ\cA(h,\one_A)
$$
and especially, we get as well an invertible $2$-cell
$\eta_A\circ f\isom\lambda_{GFA,A}(f)$; summing up, we get
an invertible $2$-cell $\eta_A\circ f\isom\one_{GFA}$. On the
other hand, the invertible modification
$(\eps*F)\odot(F*\eta)\leadsto\one_F$ provided by theorem
\ref{th_2-adjunction}(i) yields an invertible $2$-cell
$\eps_{FA}\circ F\eta_A\isom\one_{FA}$. We conclude that
$F\eta_A$ is an equivalence, and then the same holds for
$\eta_A$, by virtue of lemma \ref{lem_equiv-is-preserved}(iii).
Combining with theorem \ref{th_pseudo-nat-equiv}, we see that
(b) holds for $\eta$.

Next, obviously (b)$\Rightarrow$(c). Suppose then that
$\eta':\one_\cA\isom GF$ is a pseudo-natural equivalence,
and denote also by $i_{\cA^o}$ the identity automorphism of
$\one_{\cA^o}$; we choose a pseudo-natural equivalence
$\psi:H_\cA(\one_\cA,G)\isom H_\cB(F,\one_\cB)$ and we set
$$
\xi:=(\psi*(\one_{\cA^o}\times F))\odot(H_\cA*(i_{\cA^o}\times\eta')):
H_\cA\isom H_\cB(F,F).
$$
Since $\xi$ is pseudo-natural, for every $1$-cell $f:A'\to A$
in $\cA$, the coherence constraint $\tau^\xi$ of $\xi$ yields
invertible $2$-cells :
\set\begin{equation}\label{eq_A-and-A-prime}
Ff\circ\xi_{A'A'}(\one_{A'})\xrightarrow{\tau^\xi_{(f,\one_A),\one_{A'}}}
\xi_{A'A}(f)\xleftarrow{\tau^\xi_{(\one_{A'},f),\one_A}}
\xi_{AA}(\one_A)\circ Ff.
\end{equation}
Moreover, since $\xi_{AA}:\cA(A,A)\isom\cB(FA,FA)$ is an
equivalence for every $A\in\Ob(\cA)$, there exists a $1$-cell
$h:A\to A$ with an isomorphism $\xi_{AA}(h)\isom\one_A$. Taking
$A'=A$ and $f:=h$ in \eqref{eq_A-and-A-prime}, we conclude that
the $1$-cell $g_A:=\xi_{AA}(\one_A)$ is an equivalence in $\cA$
for every $A\in\Ob(\cA)$. On the other hand, the naturality of
$\tau^\xi$ easily implies that the rule :
$$
(f:A'\to A)\mapsto\tau^\xi_{(f,\one_A),\one_A}:g_A\circ Ff\isom\xi_{A'A}(f)
$$
defines an isomorphism of functors
$\cA(\one_A,g_A)\isom\xi_{A'A}$. But since $g_A$ is an equivalence,
$\cA(\one_A,g_A)$ is an equivalence of categories, and the same
holds for $\xi_{A'A}$. It follows that the rule : $f\mapsto Ff$
yields an equivalence $F_{A'A}:\cA(A',A)\isom\cB(FA',FA)$, {\em i.e.}
(a) holds for $F$.

Lastly, for a given $\eta'$ as in (c), there exists a pseudo-natural
equivalence $\omega:\one_\cA\isom\one_\cA$ such that
$\eta\odot\omega=\eta'$ in $\overline{2\tdu\bCat}$. We set
$\eps':=\eps\odot(F*\omega^{-1}*G):FG\Rightarrow\one_\cB$. We
compute in $\overline{2\tdu\bCat}$ :
$$
\begin{aligned}
(G*\eps')\odot(\eta'*G)&\,=
(G*\eps)\odot(GF*\omega^{-1}*G)\odot(\eta'*G) \\
&\,=(G*\eps)\odot(GF*\omega^{-1}*G)\odot(\eta*G)\odot(\omega*G) \\
&\,=(G*\eps)\odot(\eta*G)\odot(\omega^{-1}*G)\odot(\omega*G) \\
&\,=(G*\eps)\odot(\eta*G) \\
&\,=\one_G.
\end{aligned}  
$$
Likewise we get $(\eps'*F)\odot(F*\eta')=\one_F$ in
$\overline{2\tdu\bCat}$ : the details shall be left to the reader.
\end{proof}

\begin{remark}\label{rem_strict-strong}
Let $F:\cA\to\cB$ be a strict and strong $2$-equivalence between
two $2$-categories. Then there exists a strict and strong
$2$-equivalence $G:\cB\to\cA$ with strict pseudo-natural isomorphisms
$\eta:\one_\cB\isom F\circ G$ and $\eps:G\circ F\isom\one_\cA$. Indeed,
for every $B\in\Ob(\cB)$ let us pick $GB\in\Ob(\cA)$ with an isomorphism
$\eta_B:B\isom FGB$; for every $B,B'\in\Ob(\cB)$ we then get an
isomorphism of categories
$$
\cA(GB,GB')\xrightarrow{F_{GB,GB'}}\cB(FGB,FGB')
\xrightarrow{\cB(\eta_B,\eta^{-1}_{B'})}\cB(B,B')
$$
where $\cB(\eta_B,\eta^{-1}_{B'})$ denotes the functor given by
the rule : $\phi\mapsto\eta^{-1}_{B'}\circ\phi\circ\eta_B$ for
every $1$-cell $\phi:FGB\to FGB'$ of $\cB$, and
$\beta\mapsto\eta^{-1}_{B'}*\beta*\eta_B$ for every $2$-cell
$\beta$ between $1$-cells $\phi,\phi':FGB\to FGB'$. We let
$G_{B,B'}:\cB(B,B')\to\cA(GB,GB')$ be the inverse of this
isomorphism. Then the sought $G$ is given by the rules :
$B\mapsto GB$ for every $B\in\Ob(\cB)$, and
$(B,B')\mapsto G_{B,B'}$ for every pair $(B,B')$ of objects
of $\cB$ : the straightforward verification is left to the
reader. It is also clear that the rule : $B\mapsto\eta_B$
for every $B\in\Ob(\cB)$ defines a strict pseudo-natural
isomorphism $\eta$. Lastly, we define $\eps$ by the rule :
$A\mapsto F^{-1}_{GFA,A}(\eta_{FA}^{-1})$ for every $A\in\Ob(\cA)$.
The following corollary extends this simple observation to every
$2$-equivalence.
\end{remark}

\begin{corollary}\label{cor_characterize-2-eq}
Let $\cA,\cB$ be two $2$-categories, $F:\cA\to\cB$ a
pseudo-functor. We have :
\begin{enumerate}
\item
The following conditions are equivalent :
\begin{enumerate}
\alphaenu
\item
$F$ is a $2$-equivalence.
\item
There exists a pseudo-functor $L:\cB\to\cA$ with two
pseudo-natural equivalences
$$
\eta^*:\one_\cB\Rightarrow F\circ L
\qquad\text{and}\qquad
\eps^*:L\circ F\Rightarrow\one_\cA.
$$
\end{enumerate}
We say that $L$ is a {\em pseudo-inverse} for $F$.
\item
Suppose that $F$ is a $2$-equivalence. Then, for every $L$
and $\eta^*$ as in {\em (i)} there exist a pseudo-natural
equivalence $\eps^*:L\circ F\Rightarrow\one_\cA$ and
invertible modifications
$$
\Theta:(F*\eps^*)\odot(\eta^**F)\leadsto\one_F
\qquad
\Theta':(\eps^**L)\odot(L*\eta^*)\leadsto\one_L.
$$
\end{enumerate}
\end{corollary}
\begin{proof}(i.b)$\Rightarrow$(i.a): The coherence constraint
of $\eps^*$ gives an isomorphism of functors
$$
H_\cA(\one_{LFA},\eps^*_{A'})\circ L_{FA,FA'}\circ F_{AA'}\Rightarrow
H_\cA(\eps^*_A,\one_{A'})
\qquad
\text{for every $A,A'\in\Ob(\cA)$}.
$$
However, it is easily seen that both $H_\cA(\one_{LFA},\eps_{A'})$
and $H_\cA(\eps_A,\one_{A'})$ are equivalences of categories
(cp. the proof of lemma \ref{lem_long-time}), hence the same
holds for $L_{FA,FA'}\circ F_{AA'}$. Arguing likewise with the
coherence constraint for $\eta^*$, we see as well that
$F_{LB,LB'}\circ L_{BB'}$ is an equivalence for every $B,B'\in\Ob(\cB)$.
It follows easily that $F_{AA'}$ and $L_{BB'}$ are faithful for
every such $A,A'$ and $B,B'$; especially, $L_{FA,FA'}$ and
$F_{LB,LB'}$ are both faithful, and therefore $F_{AA'}$ and
$L_{BB'}$ must also be full. Lastly, let $f:FA\to FA'$ be
any $1$-cell of $\cB$; by the foregoing there exists a
$1$-cell $g:A\to A'$ with an invertible $2$-cell
$b:LFg\Rightarrow Lf$, and since $L_{FA,FA'}$ is fully
faithful, we have $b=Lc$ for an invertible $2$-cell
$c:Fg\Rightarrow f$, which shows that $F_{AA'}$ is essentially
surjective, so it is an equivalence of categories. Furthermore,
theorem \ref{th_pseudo-nat-equiv} says that the $1$-cell
$\eta^*_B:B\to FLB$ is an equivalence for every $B\in\Ob(\cB)$,
so (i.a) holds.

(i.b)$\Rightarrow$(i.a): Without loss of generality, we may
assume that both $\cA$ and $\cB$ are small $2$-categories.
Now, consider the category $\cA^F$ and the pseudo-functors
$F^\flat:\cA^F\to\cB$, $\pi^F:\cA^F\to\cA$ and $\sigma^F:\cA\to\cA^F$
provided by proposition \ref{prop_strictification}; we remark
that $\pi^F$ is a $2$-equivalence, and $\sigma^F$ is its
pseudo-inverse : indeed, we have $\pi^F\circ\sigma^F=\one_\cA$,
and a pseudo-natural equivalence
$\mu^F:\sigma^F\circ\pi^F\isom\one_\cA$, so the assertion follows
from the foregoing. It is then clear that $F$ is a $2$-equivalence
if and only if the same holds for $F\circ\pi^F$, if and only if
the same holds for $F^\flat$. Now, suppose we have found a
pseudo-functor $L':\cB\to\cA^F$, and pseudo-natural
equivalences $\eta^{**}:\one_\cB\isom F^\flat\circ L'$,
$\eps^{**}:L'\circ F^\flat\isom\one_{\cA^F}$; then set
$L:=\pi^F\circ L'$, and notice that we get pseudo-natural
equivalences
$$
(F^\flat*\mu^F*L')^{-1}\odot\eta^{**}:\one_\cB\isom F\circ L
\qquad
\pi^F*\eps^{**}*\sigma^F:L\circ F\isom\one_\cA.
$$
We may thus replace $F$ by $F^\flat$, and assume from start
that $F$ is strict. Next, consider the pseudo-natural
transformation $H_F:H_\cA\Rightarrow H_\cB(F,F)$ assigned
to $\one_F:F\Rightarrow F$, as in lemma
\ref{lem_shouldbe-towards-2-yoneda}(i).
Notice that our assumption on $F$ means that $H_{F,(A,A')}$
is an equivalence in $\bCat$, for every object $(A,A')$ of
$\cA^o\times\cA$, so theorem \ref{th_pseudo-nat-equiv} implies
that $H_F$ is a pseudo-natural equivalence of pseudo-functors.
We may therefore find a pseudo-natural equivalence
$G:H_\cB(F,F)\Rightarrow H_\cA$, with coherence constraint
$\tau^G$, and an adjunction :
$$
\eta:\one_{H_\cA}\leadsto G\odot H_F
\qquad
\eps:H_F\odot G\leadsto\one_{H_\cB\circ(F,F)}
$$
for the pair $(H_F,G)$. Explicitly, $G$ attaches to every
$A,A'\in\Ob(\cA)$ a functor
$$
G_{AA'}:\cB(FA,FA')\to\cA(A,A')
$$
which is an equivalence of categories, and $\eta$ and $\eps$
are given by natural isomorphisms of functors
$$
\eta_{AA'}:\one_{\cA(A,A')}\Rightarrow G_{AA'}\circ F_{AA'}
\qquad
\eps_{AA'}:F_{AA'}\circ G_{AA'}\Rightarrow\one_{\cB(FA,FA')}.
$$
On the other hand, by assumption, for every $B\in\Ob(\cB)$ we
may find $LB\in\Ob(\cA)$ with an equivalence $f_B:B\to FLB$;
then lemma \ref{lem_long-time}(i) says that we may also find a
$1$-cell $g_B:FLB\to B$ and invertible $2$-cells
$$
\eta'_B:\one_{FLB}\Rightarrow f_B\circ g_B
\qquad
\eps'_B:g_B\circ f_B\Rightarrow\one_B
$$
that form an adjunction for the pair $(g_B,f_B)$. To ease
notation, we set $T_{BB'}:=\cB(g_B,f_{B'})$ for every
$B,B'\in\Ob(\cB)$ and we define
$$
L_{BB'}:=G_{LB,LB'}\circ T_{BB'}:\cB(B,B')\to\cA(LB,LB').
$$
We shall henceforth drop the subscripts from the notation
for $F,G,L,T,\eta$ and $\eps$. We claim that $L$ is a
pseudo-functor $\cB\to\cA$, with coherence constraint given
by the system of $2$-cells in $\cA$
$$
\begin{aligned}
\gamma^L_{h,k}:=\,&
L(k*\eps'_{B'}*h)\odot G(\eps_{Tk}*\eps_{Th})\odot\eta_{Lk\circ Lh}:
Lk\circ Lh\Rightarrow L(k\circ h) \\
\delta^L_B:=\,&
G(\eta'_B)\odot\eta_{\one_{LB}}:\one_{LB}\Rightarrow L(\one_B)
\end{aligned}
$$
for every $B,B',B''\in\Ob(\cB)$ and every pair of $1$-cells
$h:B\to B'$ and $k:B'\to B''$ in $\cB$. Indeed, the naturality
of the rule $(h,k)\mapsto\gamma^L_{h,k}$ is clear from the
definition; to check the composition axiom, we consider
another $1$-cell $l:B''\to B'''$ in $\cB$, and we need to
show that
$$
\gamma^L_{k\circ h,l}\odot(Ll*\gamma^L_{h,k})=
\gamma^L_{h,l\circ k}\odot(\gamma^L_{l,k}*Lh).
$$
However, the naturality of $\eta$ yields the identities :
$$
\begin{aligned}
\eta_{Ll\circ L(k\circ h)}\odot(Ll*\gamma^L_{h,k})=\,&
GF(Ll*\gamma^L_{h,k})\odot\eta_{Ll\circ Lk\circ Lh} \\
\eta_{L(l\circ k)\circ Lh}\odot(\gamma^L_{k,l}*Lh)=\,&
GF(\gamma^L_{k,l}*Lh)\odot\eta_{Ll\circ Lk\circ Lh}.
\end{aligned}
$$
To ease notation, set
$$
X:=(\eps_{T(l\circ k)}*\eps_{Th})\odot F(\gamma^L_{k,l}*Lh)
\qquad
Y:=
(\eps_{Tl}*\eps_{T(k\circ h)})\odot F(Ll*\gamma^L_{h,k}).
$$
We may compute :
$$
\begin{aligned}
X=&\,
(\eps_{T(l\circ k)}*\eps_{Th})\odot
(F(L(l*\eps'_{B''}*k)\odot G(\eps_{Tl}*\eps_{Tk})\odot
\eta_{Ll\circ Lk})*FLh) \\
=&\,
(\eps_{T(l\circ k)}*\eps_{Th})\odot
(F(L(l*\eps'_{B''}*k)\odot G(\eps_{Tl}*\eps_{Tk}))*FLh)
\odot F(\eta_{Ll\circ Lk}*Lh) \\
=&\,
(\eps_{T(l\circ k)}*\eps_{Th})\odot
(FG(T(l*\eps'_{B''}*k)\odot(\eps_{Tl}*\eps_{Tk}))*FLh)
\odot F(\eta_{Ll\circ Lk}*Lh) \\
=&\,
((T(l*\eps'_{B''}*k)\odot(\eps_{Tl}*\eps_{Tk}))*Th)\odot
(\eps_{F(Ll\circ Lk)}*\eps_{Th})\odot F(\eta_{Ll\circ Lk}*Lh) \\
=&\,
((T(l*\eps'_{B''}*k)\odot(\eps_{Tl}*\eps_{Tk}))*Th)\odot
(F(Ll\circ Lk)*\eps_{Th}) \\
=&\,
((T(l*\eps'_{B''}*k)\odot(\eps_{Tl}*\eps_{Tk}))*Th)\odot
(FLl*FLk*\eps_{Th}) \\
=&\,
(T(l*\eps'_{B''}*k)*Th)\odot(\eps_{Tl}*\eps_{Tk}*\eps_{Th})
\end{aligned}
$$
where the fourth equality follows from the naturality
of $\eps$, the fifth follows from the triangular
identities \eqref{subsec_adj-pair} for the adjunction
$(\eta,\eps)$, and the sixth and seventh follow from remark
\ref{rem_equiv-2-cat}(i). Likewise, a similar calculation
gives :
$$
Y=(Tl*T(k*\eps'_{B'}*h))\odot(\eps_{Tl}*\eps_{Tk}*\eps_{Th}).
$$
To conclude, it suffices to notice the identity :
$$
T(l*k*\eps'_{B'}*h)\odot(T(l*\eps'_{B''}*k)*Th)=
T(l*\eps'_{B''}*k*h)\odot(Tl*T(k*\eps'_{B'}*h)).
$$
Next, let us show the unit axiom : we come down to checking
the identity
$$
\gamma^L_{\one_B,f}\odot(Lf*\delta^L_B)=\one_{Lf}
\qquad 
\text{for every $1$-cell $f:B\to B'$ in $\cB$}.
$$
However,
$$
\begin{aligned}
\gamma^L_{\one_B,f}\odot(Lf*\delta^L_B)=\,&
L(f*\eps'_B)\odot G(\eps_{Tf}*\eps_{T\one_B})\odot
\eta_{Lf\circ L\one_B}\odot(Lf*\delta^L_B) \\
=\,&
L(f*\eps'_B)\odot G(\eps_{Tf}*\eps_{T\one_B})\odot
GF(Lf*\delta^L_B)\odot\eta_{Lf}
\end{aligned}
$$
by the naturality of $\eta$; taking into account the
triangularity identities \eqref{subsec_adj-pair} for
the adjunction $(\eta,\eps)$, we are then reduced to
showing that
$$
Z:=T(f*\eps'_B)\odot(\eps_{Tf}*\eps_{T\one_B})\odot
F(Lf*\delta^L_B)=\eps_{Tf}.
$$
However, from the triangular identities for the adjunction
$(\eta'_B,\eps'_B)$, we easily obtain :
$$
T(f*\eps'_B)=Tf*\eta'^{-1}_B
$$
and on the other hand, the naturality of $\eps$ yields
$$
\eta'^{-1}_B\odot\eps_{T\one_B}=\eps_{\one_{FLB}}\odot FG\eta'^{-1}_B
$$
whence :
$$
\begin{aligned}
Z=\,&
(\eps_{Tf}*(\eps_{\one_{FLB}}\odot FG\eta'^{-1}_B))\odot
(FLf*F\delta^L_B) \\
=\,&
(\eps_{Tf}*(\eps_{\one_{FLB}}\odot FG\eta'^{-1}_B))\odot
(FLf*FG\eta'_B)\odot(FLf*F\eta_{\one_{LB}}) \\
=\,&
(\eps_{Tf}*\eps_{\one_{FLB}})\odot(FLf*F\eta_{\one_{LB}}).
\end{aligned}
$$
From the triangular identities for $(\eta,\eps)$ we
also see that
$F(\eta_{\one_{LB}})=\eps^{-1}_{F\one_{LB}}=\eps^{-1}_{\one_{FLB}}$,
so finally :
$$
Z=(\eps_{Tf}*\eps_{\one_{FLB}})\odot(FLf*\eps^{-1}_{\one_{FLB}})
=\eps_{Tf}
$$
as required. A similar computation establishes the second
identity for the unit axiom of $L$. This completes the
construction of the pseudo-functor $L$. Next we set, for
every $B,B'\in\Ob(\cB)$ and every $1$-cell $h:B\to B'$
in $\cB$
$$
\begin{aligned}
\eta^*_B:=&\, f_B:B\to FLB \\
\tau^*_h:=&\, (f_{B'}*h*\eps'_B)\odot(\eps_{Th}*f_B):
FL(h)\circ f_B\Rightarrow f_{B'}\circ h.
\end{aligned}
$$

\begin{claim} The system of $1$-cells $\eta^*_\bullet$
defines a pseudo-natural equivalence
$\eta^*:\one_\cB\Rightarrow FL$ with coherence constraint
given by the system of $2$-cells $\tau^*_\bullet$.
\end{claim}
\begin{pfclaim} By remark \ref{rem_pseudo-natural}(i), the
naturality of $\tau^*$ amounts to the identity
$$
(f_{B'}*\beta)\odot(f_{B'}*h_1*\eps'_B)\odot(\eps_{Th_1}*f_B)
=(f_{B'}*h_2*\eps'_B)\odot(\eps_{Th_2}*f_B)\odot(FL(\beta)*f_B)
$$
for every pair of $1$-cells $h_1,h_2:B\to B'$ and every
$2$-cell $\beta:h_1\Rightarrow h_2$. However, a simple
inspection shows that
$$
(f_{B'}*\beta)\odot(f_{B'}*h_1*\eps'_B)=
(f_{B'}*h_2*\eps'_B)\odot(T(\beta)*f_B)
$$
so we are reduced to checking that
$$
T(\beta)\odot\eps_{Th_1}=\eps_{Th_2}\odot FL\beta
$$
which in turns follows from the naturality of $\eps$.
Next, let us check the coherence axioms for $\tau^*$.
Denote by $(\delta^{FL},\gamma^{FL})$ the coherence
constraint of $FL$ (see remark \ref{rem_pseudo-funct}(v));
for given $1$-cells $h:B\to B'$ and $h':B'\to B''$ of
$\cB$ we have to verify the identities
\set\begin{equation}\label{eq_tau-star}
\tau^*_{\one_B}\odot(\delta^{FL}_B*f_B)=\one_{f_B}
\qquad
(\tau^*_{h'}*h)\odot(FLh'*\tau^*_h)
=\tau^*_{h'\circ h}\odot(\gamma^{FL}_{h,h'}*f_B).
\end{equation}
However, a direct inspection gives the identity
$$
(\tau^*_{h'}*h)\odot(FLh'*\tau^*_h)=
(f_{B''}*h'*\eps'_{B'}*h*\eps'_B)\odot(\eps_{Th'}*\eps_{Th}*f_B)
$$
so we are reduced to showing that
$$
T(h'*\eps'_{B'}*h)\odot(\eps_{Th'}*\eps_{Th})
=\eps_{T(h'\circ h)}\odot\gamma^{FL}_{h,h'}.
$$
Now, we compute :
$$
\begin{aligned}
\eps_{T(h'\circ h)}\odot\gamma^{FL}_{h,h'}
=&\,
\eps_{T(h'\circ h)}\odot F(\gamma^L_{h,h'}) \\
=&\,
\eps_{T(h'\circ h)}\odot FG(T(h'*\eps'_{B'}*h)\odot(\eps_{Th'}*\eps_{Th}))
\odot F(\eta_{Lh'\circ Lh}) \\
=&\,
T(h'*\eps'_{B'}*h)\odot(\eps_{Th'}*\eps_{Th})\odot
\eps_{F(Lh'\circ Lh)}\odot F(\eta_{Lh'\circ Lh}) \\
=&\,T(h'*\eps'_{B'}*h)\odot(\eps_{Th'}*\eps_{Th})
\end{aligned}
$$
where the third identity follows from the naturality of $\eps$,
and the fourth follows from the triangular identities
\eqref{subsec_adj-pair} for the pair $(\eta,\eps)$. For
the first identity \eqref{eq_tau-star}, we notice that
$$
\tau^*_{\one_B}=(f_B*\eps'_B)\odot(\eps_{T\one_B}*f_B)=
(\eta'^{-1}_B*f_B)\odot(\eps_{T\one_B}*f_B)
$$
due to the triangular identities for the pair $(\eta',\eps')$,
so we are reduced to checking that
$$
\eps_{T\one_B}\odot\delta^{FL}_B=\eta'_B.
$$
We compute :
$$
\eps_{T\one_B}\odot\delta^{FL}_B=\eps_{T\one_B}\odot F\delta^L_B=
\eps_{T\one_B}\odot F(G(\eta'_B)\odot\eta_{\one_{LB}})=
\eta'_B\odot\eps_{\one_{FLB}}\odot F(\eta_{\one_{LB}})=
\eta'_B
$$
which concludes the verification of the pseudo-naturality
of $\eta^*$. Lastly, $\eta^*_B$ is an equivalence for every
$B\in\Ob(\cB)$, so $\eta^*$ is a pseudo-natural equivalence,
by theorem \ref{th_pseudo-nat-equiv}.
\end{pfclaim}

\begin{claim}\label{cl_G-is-a-2-equiv}
For every $A\in\Ob(\cA)$, the $1$-cell $G(g_{FA}):LFA\to A$
is an equivalence in $\cA$.
\end{claim}
\begin{pfclaim} We have an invertible $2$-cell
$\eps_{g_{FA}}:FG(g_{FA})\Rightarrow g_{FA}$, and since
$g_{FA}$ is an equivalence, it follows that the same holds
for $FG(g_{FA})$. Then the assertion follows from lemma
\ref{lem_equiv-is-preserved}(iii).
\end{pfclaim}

A simple inspection shows that $L_{BB'}$ is equivalence
of categories, for every $B,B'\in\Ob(\cB)$. Taking into
account claim \ref{cl_G-is-a-2-equiv}, we deduce that $L$
is a $2$-equivalence; therefore, the foregoing discussion
applies as well with $F$ replaced by $L$, and yields
another $2$-equivalence $M:\cA\to\cB$ with a pseudo-natural
equivalence $\eta^{**}:\one_\cA\Rightarrow LM$. Pick also
pseudo-natural equivalences
$$
\mu^*:FL\Rightarrow\one_\cB
\qquad\text{and}\qquad
\mu^{**}:LM\Rightarrow\one_\cA
$$
representing the inverses of the $2$-cells $\eta^*$, $\eta^{**}$
in $\overline{2\tdu\bCat}$. We get a pseudo-natural equivalence
$$
\beta:=(\mu^**M)\odot(F*\eta^{**}):F\Rightarrow M
$$
(corollary \ref{cor_first-corollary}), whence the sought
pseudo-natural equivalence
$$
\eps^*:=\mu^{**}\odot(L*\beta):LF\Rightarrow\one_\cA.
$$

(ii): Assertion (i) says that $F$ is a $2$-equivalence if and
only if it is an equivalence in the $2$-category
$\overline{2\tdu\bCat}$. Then the assertion follows immediately
from lemma \ref{lem_long-time}(i).
\end{proof}

\subsection{\texorpdfstring{$2$}{2}-Limits and
\texorpdfstring{$2$}{2}-colimits}\label{sec_pseudo-cones}
Consider $2$-categories $\cA$, $\cB$. For every object
$B$ of $\cB$, one may define the
{\em constant pseudo-functor\/} with value $B$ : this
is the pseudo-functor
$$
\sF_B:\cA\to\cB
\qquad\text{such that}\qquad
\sF_B(A):=B
\qquad
\sF_B(f):=\one_B
\qquad
\sF_B(\alpha):=i_B
$$
for every $A\in\Ob(\cA)$, every $1$-cell $f$, and every
$2$-cell $\alpha$ of $\cA$. The coherence constraint for
$\sF_B$ consists of identities. Given a pseudo-functor
$F:\cA\to\cB$, a {\em pseudo-cone on $F$ with vertex $B$\/}
is a pseudo-natural transformation $\sF_B\Rightarrow F$.
Dually, a {\em pseudo-cocone on $F$ with vertex $B$\/}
is a pseudo-natural transformation $F\Rightarrow\sF_B$.
Especially, every $1$-cell $f:B\to B'$ induces a pseudo-cone
with vertex $B$ :
$$
\sF_f:\sF_B\Rightarrow\sF_{B'}
\quad :\quad
(\sF_f)_A:=f
\qquad
\text{for every $A\in\Ob(\cA)$}
$$
whose coherence constraint consists of identities. Notice
that $\sF_f$ can also be viewed as a pseudo-cocone with
vertex $B'$. Also, every $2$-cell $\beta:f\Rightarrow f'$
induces a modification
$$
\sF_\beta:\sF_f\leadsto\sF_{f'}
\quad :\quad
A\mapsto\beta
\qquad
\text{for every $A\in\Ob(\cA)$}.
$$

\begin{definition}\label{def_pseudo-lim}
Let $I$ and $\cB$ be two $2$-categories, $F:I\to\cB$ a pseudo-functor.

(i)\ \
A {\em $2$-limit\/} (resp. a {\em strong $2$-limit}) of $F$ is a pair :
$$
\Pslim{I}F:=(L,\pi)
$$
consisting of an object $L$ of $\cB$ and a pseudo-cone
$\pi:\sF_L\Rightarrow F$, such that the functor :
$$
\cB(B,L)\to\sPsNat(\sF_B,F)
\quad :\quad
f\mapsto\pi\odot\sF_f
\qquad
(\beta:f\Rightarrow f')\mapsto(\pi*\sF_\beta:\sF_f\leadsto\sF_{f'})
$$
is an equivalence (resp. an isomorphism) of categories, for every
$B\in\Ob(\cB)$. In this case, we also say that $\pi$ is a
{\em universal pseudo-cone} (resp. a {\em strong universal
pseudo-cone}).

(ii)\ \
Dually, a {\em $2$-colimit\/} (resp. a {\em strong $2$-colimit})
of $F$ is a pair :
$$
\Pscolim{I}F:=(L,\pi)
$$
consisting of an object $L$ of $\cB$ and a pseudo-cocone
$\pi:F\Rightarrow\sF_L$, such that the functor:
$$
\cB(L,B)\to\sPsNat(F,\sF_B)
\quad :\quad
f\mapsto\sF_f\odot\pi
\qquad
(\beta:f\Rightarrow f')\mapsto(\sF_\beta*\pi:\sF_f\leadsto\sF_{f'})
$$
is an equivalence (resp. an isomorphism) of categories, for every
$B\in\Ob(\cB)$. In this case, we also say that $\pi$ is a
{\em universal pseudo-cocone} (resp. a {\em strong universal
pseudo-cocone}).

(iii)\ \
We say that $\cB$ is {\em $2$-complete} (resp.
{\em $2$-cocomplete}, resp. {\em strongly $2$-complete},
resp. {\em strongly $2$-cocomplete}) if, for every small
$2$-category $I$, every pseudo-functor $I\to\cB$ admits
a $2$-limit (resp. a $2$-colimit, resp. a strong $2$-limit,
resp. a strong $2$-colimit).

(iv)\ \
Let $\cC$ be any other $2$-category, $G:\cB\to\cC$ any
pseudo-functor, and $(L,\pi)$ a $2$-limit (resp. a
$2$-colimit) for $F$. We say that {\em $G$ commutes
with the $2$-limit of $F$} (resp. that {\em $G$ commutes
with the $2$-colimit of $F$}) if the pseudo-cone
$G*\pi:\sF_{GL}\Rightarrow G\circ F$ (resp. the pseudo-cocone
$G*\pi:G\circ F\Rightarrow\sF_{GL}$) is universal.
\end{definition}

\begin{remark}\label{rem_pseudo-limit}
(i)\ \
As usual, if the $2$-limit of a pseudo-functor exists, it is
unique up to (non-unique) equivalence : if $(L',\pi')$ is
another $2$-limit, there exist an equivalence $h:L\to L'$
and an isomorphism $\beta:\pi'\odot\sF_h\isom\pi$; moreover,
the pair $(h,\beta)$ is unique up to unique isomorphism, in
a suitable sense, that the reader may spell out, as an
exercise. A similar remark holds for $2$-colimits.

(ii)\ \
Likewise, in the situation of definition
\ref{def_pseudo-lim}(iv), the commutation of $G$ with the
$2$-limit or $2$-colimit of $F$ is an intrinsic property
of $G$, {\em i.e.}  is independent of the choice of a
$2$-limit or $2$-colimit for $F$.

(iii)\ \
To be in keeping with the terminology of \cite[Ch.VII]{Bor},
we should write pseudo-bilimit instead of 2-limit (and
likewise for 2-colimit). The term ``$2$-limit'' denotes
in {\em loc.cit.} a related notion, which is unique
up to isomorphism, not just up to equivalence. However,
the notion introduced in definition \ref{def_pseudo-lim}
occurs more frequently in applications.

(iv)\ \
Let $A_\bullet:=(A_i~|~i\in I)$ be any family of objects of
the $2$-category $\cA$, indexed by a small set $I$. We may
regard $A_\bullet$ as a strict pseudo-functor from the discrete
category $I$ to $\cA$. Then a $2$-limit (resp. a $2$-colimit)
of this pseudo-functor shall also be called a {\em $2$-product}
(resp. a $2$-coproduct) of $A_\bullet$.

(v)\ \
Let $L$ be the small category with $\Ob(L):=\{0,1,2\}$,
and whose set of arrows consists of the identity morphisms,
and two more arrows $1\to 0$ and $2\to 0$.
An essentially commutative square \eqref{eq_2-cell}
can be regarded as a pseudo-cone $\pi$ with vertex $A$,
on the functor $F:L\to\cA$ such that $F(0):=D$, $F(1):=B$,
$F(2):=C$, $F(1\to 0):=h$ and $F(2\to 0):=k$.
We say that \eqref{eq_2-cell} is {\em $2$-cartesian}
if $(A,\pi)$ is a $2$-limit of the functor $F$. This
$2$-limit shall be called the {\em $2$-fibre product}
of $k$ and $h$, and shall be denoted
$$
B\mathop{\times}^2_{(h,k)}C
$$
or sometimes, just $B\stackrel{2}{\times}_DC$, if there
is no danger of ambiguity.

(vi)\ \
In view of \eqref{subsec_opposite-mods}, it is easily seen
that a pair $(L,\pi)$ is a $2$-colimit for the pseudo-functor
$F:I\to\cB$ if and only if the pair $(L^o,\pi^o)$ is a
$2$-limit for the pseudo-functor $F^o:I^o\to\cB^o$, and if and
only if $({}^oL,{}^o\pi)$ is a $2$-limit for the pseudo-functor
${}^oF:{}^oI\to{}^o\cB$. Especially, $\cB$ is $2$-cocomplete if
and only if $\cB^o$ is $2$-complete, and if and only if ${}^o\cB$
is $2$-cocomplete.

(vii)\ \
Let $I$ and $\cB$ be two categories, which we regard as
$2$-categories with trivial $2$-cells. Then every pseudo-functor
$F:I\to\cB$ is a functor, every pseudo-natural transformation
$F\Rightarrow G$ of functors $F,G:I\to\cB$ is a natural
transformation, and the category $\sPsNat(F,G)$ is discrete.
It follows easily that a $2$-limit $(L,\pi)$ of any such
functor $F:I\to\cB$ is a strong $2$-limit, it represents
also the (usual) limit of $F$, and the universal pseudo-cone
$\pi$ is also a universal cone. Likewise, any $2$-colimit of
$F$ represents the colimit of $F$, and every universal
pseudo-cocone is a universal cocone. 
\end{remark}

The following lemma \ref{lem_pseudo-trivial} indicates that
the framework of $2$-categories does indeed provide an adequate
answer to the issues raised in the introduction of this chapter.

\begin{lemma}\label{lem_pseudo-trivial}
Let $I$ be a $2$-category, $F,G:I\to\cB$ two
pseudo-functors, $\omega:F\Rightarrow G$ a pseudo-natural
equivalence, and suppose that the $2$-limit of $F$ exists.
Then the same holds for the $2$-limit of\/ $G$, and there is
a natural equivalence in $\cB$ :
$$
\Pslim{I}F\isom\Pslim{I}G.
$$
More precisely, if $(L,\pi)$ is a pair with $L\in\Ob(\cB)$
and a pseudo-cone $\pi:\sF_L\Rightarrow F$ representing the
$2$-limit of $F$, then the pair $(L,\omega\odot\pi)$
represents the $2$-limit of $G$.

A dual assertion holds for $2$-colimits.
\end{lemma}
\begin{proof} It is easily seen that the rule
$\alpha\mapsto\omega\odot\alpha$ induces an equivalence
of categories :
$$
\sPsNat(\sF_B,F)\to\sPsNat(\sF_B,G)
\qquad
\text{for every $B\in\Ob(\cB)$}.
$$
The claim is an immediate consequence.
\end{proof}

\sset\subsubsection{}\label{subsec_Yoneda-prsrvs-univ-cone}
Let $\cC$ be any $2$-category with small $\Hom$-categories.
To every object $X$ of $\cC$ we attach the strict pseudo-functor
$$
h_X:\cC^o\to\bCat
\qquad
Y\mapsto\cC(Y,X)\subset Y/\cC
$$
that assigns to every $1$-cell $Y'\xrightarrow{f}Y$ of $\cC$
the restriction $\cC(Y,X)\to\cC(Y',X)$ of the strict
pseudo-functor $f^*:Y/\cC\to Y'/\cC$, and to every $2$-cell
$\beta:f\Rightarrow f'$ the restriction of the strict
pseudo-natural transformation $\beta^*:f^*\Rightarrow f'^*$
(notation of example \ref{ex_pseudo-functors}(i)).

Notice that $h_{X^o}:\cC\to\bCat$ is the strict pseudo-functor
given by the rule : $Y\mapsto\cC(X,Y)$ for every $Y\in\Ob(\cC)$,
and which assigns to the $1$-cell $f$ as in the foregoing the
restriction of the strict pseudo-functor $f_*:\cC/Y\to\cC/Y'$,
and to the $2$-cell $\beta$ the restriction of
$\beta_*:f_*\Rightarrow f'_*$ (again, with the notation of
example \ref{ex_pseudo-functors}(i)).

\begin{proposition}\label{prop_2-lims-and-Homs}
With the notation of \eqref{subsec_Yoneda-prsrvs-univ-cone},
let $I$ be another $2$-category, $F:I\to\cC$ a pseudo-functor,
and $L\in\Ob(\cC)$ any object. We have :
\begin{enumerate}
\item
For every pseudo-cone $\pi:\sF_L\Rightarrow F$, the following
conditions are equivalent :
\begin{enumerate}
\item
The pseudo-cone $\pi$ is universal.
\item
For every $X\in\Ob(\cC)$ the pseudo-cone
$h_{X^o}*\pi:\sF_{\cC(X,L)}\Rightarrow h_{X^o}\circ F$ is universal.
\end{enumerate}
\item
For every pseudo-cocone $\pi:F\Rightarrow\sF_L$, the following
conditions are equivalent :
\begin{enumerate}
\item
The pseudo-cocone $\pi$ is universal.
\item
For every $X\in\Ob(\cC)$ the pseudo-cone
$h_X*\pi:\sF_{\cC(L,X)}\Rightarrow h_X\circ F$ is universal.
\end{enumerate}
\end{enumerate}
\end{proposition}
\begin{proof}(i.a)$\Rightarrow$(i.b): We have to check that
for every small category $\cB$, the functor
$$
T:\bFun(\cB,\cC(X,L))\to\sPsNat(\sF_\cB,h_{X^o}*\pi)
\qquad
(G:\cB\to\cC(X,L))\mapsto(h_{X^o}*\pi)\odot\sF_G
$$
is an equivalence. To this aim, consider any pseudo-cone
$\lambda:\sF_\cB\Rightarrow h_{X^o}*\pi$.
Explicitly, $\lambda$ is the datum of a system of functors
$(\lambda_i:\cB\to\cC(X,Fi)~|~i\in\Ob(I))$ and a system of
natural isomorphisms of functors
$(\tau^\lambda_\phi:h_{X^o}(F\phi)\circ\lambda_i\Rightarrow\lambda_j~|~
(i\xrightarrow{\phi}j)\in\rMorph(I))$ fulfilling the usual
coherence axioms. To every $B\in\Ob(\cB)$ we attach the system
$$
\lambda_\bullet(B):=(\lambda_i(B):X\to Fi~|~i\in\Ob((I)).
$$
It is easily seen that $\lambda_\bullet(B):\sF_X\Rightarrow F$
is a pseudo-cone with coherence constraint given by the system
$(\tau^\lambda_{\phi,B}:F(\phi)\circ\lambda_i(B)\Rightarrow\lambda_j(B)
~|~(i\xrightarrow{\phi}j)\in\rMorph(I))$. Since $(L,\pi)$ is
universal, we may then find a $1$-cell $\lambda^\dagger_B:X\to L$
and an isomorphism of pseudo-cones
$$
\Omega^\lambda_B:\pi\odot\sF_{\lambda^\dagger_B}\leadsto\lambda_\bullet(B)
\qquad
\text{for every $B\in\Ob(\cB)$}.
$$
Moreover, for every morphism $f:B\to B'$ in $\cB$, the system
$(\lambda_i(f):\lambda_i(B)\Rightarrow\lambda_i(B')~|~i\in\Ob(I))$
yields a modification
$\lambda_\bullet(f):\lambda_\bullet(B)\leadsto\lambda_\bullet(B')$.
Then, there exists a unique $2$-cell in $\cC$
\set\begin{equation}\label{eq_for-lambda-dagger}
\lambda^\dagger_f:\lambda^\dagger_B\Rightarrow\lambda^\dagger_{B'}
\qquad\text{such that}\qquad
\Omega^\lambda_{B'}\odot(\pi*\sF_{\lambda^\dagger_f})=
\lambda_\bullet(f)\odot\Omega^\lambda_B.
\end{equation}
Clearly $\lambda_\bullet(g)\odot\lambda_\bullet(f)=
\lambda_\bullet(g\circ f)$ for every composable pair of
morphisms $B\xrightarrow{f}B'\xrightarrow{g}B''$ of $\cB$,
whence $\lambda^\dagger_g\odot\lambda^\dagger_f=\lambda^\dagger_{g\circ f}$.
We have thus associated with $\lambda$ a well defined functor
$\lambda^\dagger:\cB\to\cC(X,L)$, together with an isomorphism
$$
\Omega^\lambda_\bullet:
(h_{X^o}*\pi)\odot\sF_{\lambda^\dagger}\isom\lambda
\qquad
\text{in $\sPsNat(\sF_\cB,h_{X^o}\circ F)$}.
$$
This proves that $T$ is essentially surjective. Next, let us
check that $T$ is faithful : indeed, consider two functors
$G,H:\cB\to\cC(X,L)$ and two natural transformations
$\alpha,\beta:G\Rightarrow H$ such that
$(h_{X^o}*\pi)\odot\sF_\alpha=(h_{X^o}*\pi)\odot\sF_\beta$; the
latter identity means that
$$
\pi_i*\alpha_B=\pi_i*\beta_B:\pi_i\circ GB\Rightarrow\pi_i\circ HB
\qquad
\text{for every $i\in\Ob(I)$ and every $B\in\Ob(\cB)$}.
$$
In other words, for every $B\in\Ob(\cB)$, the morphisms of
pseudo-cones
$\pi*\sF_{\alpha_B},\pi*\sF_{\beta_B}:
\pi\odot\sF_{GB}\Rightarrow\pi\odot\sF_{HB}$ coincide; then
the universality of $\pi$ implies that $\alpha_B=\beta_B$
for every $B\in\Ob(\cB)$, whence the assertion. Lastly, to
check that $T$ is full we consider, for $G$ and $H$ as in
the foregoing, any modification
$\Xi:(h_{X^o}*\pi)\odot\sF_G\leadsto(h_{X^o}*\pi)\odot\sF_H$;
hence, $\Xi$ amounts to a system of natural transformations
$(\Xi_i:\pi_{i,*}\circ G\Rightarrow\pi_{i,*}\circ H~|~i\in\Ob(I))$
related by the compatibility conditions :
$$
\Xi_j=(\tau^\pi_{\phi*}*H)\odot((F\phi)_**\Xi_i)
\qquad
\text{for every $1$-cell $\phi:i\to j$ of $I$}
$$
where $\pi_{i,*}:\cC(X,L)\to\cC(X,Fi)$ and
$(F\phi)_*:\cC(X,Fi)\to\cC(X,Fj)$ are the functors induced
by $\pi_i:L\to Fi$ and respectively $F\phi:Fi\to Fj$, and where
$\tau^\pi_{\phi*}:(F\phi)_*\circ\pi_{i*}\isom\pi_{j*}$ denotes the
isomorphism of functors induced by the coherence constraint
$\tau^\pi$ of $\pi$. Then it is easily seen that for every
$B\in\Ob(\cB)$ the system
$$
(\Xi_{i,B}:\pi_i\circ GB\Rightarrow\pi_i\circ HB~|~i\in\Ob(I))
$$
yields a morphism of pseudo-cones
$\Xi_{\bullet,B}:\pi*\sF_{GB}\leadsto\pi*\sF_{HB}$, and the naturality
of $\Xi_i$ gives the identities
\set\begin{equation}\label{eq_clueless}
(\pi*\sF_{Hf})\odot\Xi_{\bullet,B}=\Xi_{\bullet,B'}\odot(\pi*\sF_{Gf})
\qquad
\text{for every morphism $f:B\to B'$ of $\cB$}.
\end{equation}
By the universality of $\pi$, the morphisms $\Xi_{\bullet,B}$
corresponds to a unique $2$-cell $\Xi^*_B:GB\Rightarrow HB$,
and \eqref{eq_clueless} implies that $(\Xi^*_B~|~B\in\Ob(\cB))$
is a well defined natural transformation $\Xi^*:G\Rightarrow H$.
By direct inspection, we see that $(h_{X^o}*\pi)*\sF_{\Xi^*}=\Xi$,
whence the contention.

(i.b)$\Rightarrow$(i.a): To any $X\in\Ob(\cC)$ and any pseudo-cone
$\phi:\sF_X\Rightarrow F$ we attach the pseudo-cone
$$
\phi_*:\sF_{\cC(X,X)}\Rightarrow h_{X^o}\circ F
\qquad
i\mapsto(\phi_{i*}:\cC(X,X)\to\cC(X,Fi))
$$
with coherence constraint given by the system of isomorphisms
of functors
$$
(\tau^\phi_{\psi*}:(F\psi)_*\circ\phi_{i*}\Rightarrow\phi_{j*}~|~
(i\xrightarrow{\psi}j)\in\rMorph(I))
$$
where $\tau^\phi_\bullet$ denotes the coherence constraint
of $\phi$. By assumption, there exist a functor
$\phi^\dagger:\cC(X,X)\to\cC(X,L)$ and an isomorphism of
pseudo-cones
$\omega:\phi_*\isom(h_{X^o}*\pi)\odot\sF_{\phi^\dagger}$. Then
it is easily seen that the system
$$
(\omega_{i,\one_X}:\phi_i\Rightarrow
\pi_i\circ\phi^\dagger(\one_X)~|~i\in\Ob(I))
$$
yields an isomorphism of pseudo-cones
$\phi\isom\pi\odot\sF_{\phi^\dagger(\one_X)}$. This shows that
the functor
\set\begin{equation}\label{eq_tongue-in-llama}
\cC(X,L)\to\sPsNat(\sF_X,F)
\end{equation}
as in definition \ref{def_pseudo-lim}(i) is essentially surjective.
Next, in order to check that \eqref{eq_tongue-in-llama} is faithful,
consider two $1$-cells $g,h:X\to L$ in $\cC$ and two $2$-cells
$\alpha,\beta:g\Rightarrow h$ such that $\pi_i*\alpha=\pi_i*\beta$
for every $i\in\Ob(I)$; let also $\cB$ be the initial object
of $\bCat$, {\em i.e.} $\cB$ is a category with a unique object
$b$ and a unique morphism $\one_b$. We let $G,H:\cB\to\cC(X,L)$
be the functors such that $Gb:=g$ and $Hb:=h$; then $\alpha$ and
$\beta$ correspond to natural transformations
$\alpha',\beta':G\Rightarrow H$ such that
$(h_{X^o}*\pi)\odot\sF_{\alpha'}=(h_{X^o}*\pi)\odot\sF_{\beta'}$, and
the universality of $h_{X^o}*\pi$ implies that $\alpha'=\beta'$,
whence $\alpha=\beta$, as required. Lastly, for $g$, $h$, $\cB$,
$G$ and $H$ as in the foregoing, let
$\Xi:\pi\odot\sF_g\leadsto\pi\odot\sF_h$ be any modification; we
deduce a modification
$\Xi':(h_{X^o}*\pi)\odot\sF_G\leadsto(h_{X^o}*\pi)\odot\sF_H$,
whence -- by the universality of $(h_{X^o}*\pi)$ -- a unique
natural transformation $\beta':G\Rightarrow H$ such that
$\Xi'=(h_{X^o}*\pi)\odot\sF_\beta$. Then $\beta'$ is the datum of
a $2$-cell $\beta:g\Rightarrow h$ such that $\Xi=\pi*\sF_\beta$,
which shows that the functor \eqref{eq_tongue-in-llama} is also full.

Assertion (ii) is easily deduced from (i), by considering the
pseudo-cone $\pi^\circ:\sF_{L^\circ}\Rightarrow F^\circ$ (details
left to the reader).
\end{proof}

\begin{proposition}\label{prop_2-adjs-and-lims}
Let $I$, $\cB$ and $\cC$ be three $2$-categories, $G:\cB\to\cC$,
and $H:\cC\to\cB$ two pseudo-functors, and suppose that $(H,G)$
is a $2$-adjoint pair. We have :
\begin{enumerate}
\item
The pseudo-functor $G$ commutes with the $2$-limit of every
pseudo-functor $F:I\to\cB$.
\item
Dually, $H$ commutes with the $2$-colimit of every
pseudo-functor $F':I\to\cC$.
\end{enumerate}
\end{proposition}
\begin{proof}(i): Let $F:I\to\cB$ be any pseudo-functor, and
$(L,\pi)$ any $2$-limit for $F$. We fix a unit and a counit
$(\eta,\eps)$ for the $2$-adjoint pair $(H,G)$, and a pair
of invertible modifications
$$
\Sigma:(G*\eps)\odot(\eta*G)\leadsto\one_G
\qquad
\Sigma':(\eps*H)\odot(H*\eta)\leadsto\one_H.
$$
For every $Y\in\Ob(\cC)$ we consider the diagram of functors :
$$
\cD
\qquad : \qquad
{\diagram \cC(Y,GL) \ar[r]^-\theta \ar[d] & \cB(HY,L) \ar[d] \\
\sPsNat(\sF_Y,GF) \ar[r]^-{\theta^\dagger} & \sPsNat(\sF_{HY},F)
\enddiagram}
$$
where $\theta$ is given by the rules :
$$
(f:Y\to GL)\mapsto(\eps_L\circ Hf:HY\to HGL)
\qquad\text{and}\qquad
(\beta:f\Rightarrow f')\mapsto\eps_L*H\beta
$$
for every $1$-cell $f$ and every $2$-cell $\beta$ of $\cC(Y,GL)$.
The functor $\theta^\dagger$ assigns to every pseudo-natural
transformation $\psi:\sF_Y\Rightarrow GF$ the pseudo-natural
transformation $(\eps*F)\odot(H*\psi):\sF_{HY}\Rightarrow F$,
and to every modification $\Xi:\psi\leadsto\psi'$ the
modification $(\eps*F)*(H\circ\Xi)$ (notation of remark
\ref{rem_construct-modifs}(vi)). The left (resp. right) vertical
arrow is the functor associated with the pseudo-cone $G*\pi$
(resp. $\pi$) as in definition \ref{rem_construct-modifs}(i).
By inspecting the proof of theorem \ref{rem_construct-modifs}(ii),
it is easily seen that $\theta$ is an equivalence, and by
assumption, the same holds for the right vertical arrow.
We notice :
\begin{claim}\label{cl_climate-march}
Diagram $\cD$ is essentially commutative.
\end{claim}
\begin{pfclaim} The composition of $\theta$ with the
right vertical arrow assigns to every $1$-cell $f:Y\to GL$
the pseudo-natural transformation
$\pi\odot\sF_{\eps_L\circ Hf}=\pi\odot\sF_{\eps_L}\odot\sF_{Hf}=
\pi\odot(\eps*\sF_L)\odot\sF_{Hf}$.
On the other hand, the composition of $\theta^\dagger$ with
the left vertical arrow assigns to $f$ the pseudo-natural
transformation $(\eps*F)\odot(H*((G*\pi)\odot\sF_f))$.
However, we have invertible modifications:
$$
\Theta_f:(\eps*F)\odot(H*((G*\pi)\odot\sF_f))\leadsto
(\eps*F)\odot(HG*\pi)\odot(H*\sF_f)\leadsto
\pi\odot(\eps*\sF_L)\odot\sF_{Hf}
$$
for every such $f$ (example \ref{ex_modifications}(i,ii)).
We need to check that the rule : $f\mapsto\Theta_f$ yields
a natural transformation; by unwinding the definitions, we
come down to showing the commutativity of the diagram :
$$
\xymatrix{ \eps_{Fi}\circ H(G\pi_i\circ f)
\ar@{=>}[d]_{\eps_{Fi}*\gamma^H_{f,G\pi_i}}
\ar@{=>}[rrr]^-{\eps_{Fi}*H(G\pi_i*\beta)} & & &
\eps_{Fi}\circ H(G\pi_i\circ f')
\ar@{=>}[d]^{\eps_{Fi}*\gamma^H_{f',G\pi_i}} \\
\eps_{Fi}\circ HG\pi_i\circ Hf
\ar@{=>}[rrr]^-{\eps_{Fi}*HG\pi_i*H\beta}
\ar@{=>}[d]_{(\tau^\eps_{\pi_i})^{-1}*Hf} & & &
\eps_{Fi}\circ HG\pi_i\circ Hf'
\ar@{=>}[d]^{(\tau^\eps_{\pi_i})^{-1}*Hf'} \\
\pi_i\circ\eps_L\circ Hf
\ar@{=>}[rrr]^-{(\pi_i\circ\eps_L)*H\beta} & & &
\pi_i\circ\eps_L\circ Hf'
}$$
where $(\delta^H,\gamma^H)$ is the coherence constraint of
$H$ and $\tau^\eps$ is the coherence constraint of $\eps$.
However, the commutativity of the top square subdiagram
follows from the naturality of $\gamma^H$ (remark
\ref{rem_pseudo-funct}(ii)). The commutativity of the bottom
square subdiagram is obvious, since the compositions of both
pairs of arrows equal $(\tau^\eps_{\pi_i})^{-1}*H\beta$ (remark
\ref{rem_equiv-2-cat}(i)).
\end{pfclaim}

In view of claim \ref{cl_climate-march}, we are reduced to
checking that $\theta^\dagger$ is an equivalence of categories.
To this aim, let us consider the functor
$$
\mu:\sPsNat(\sF_{HY},F)\to\sPsNat(\sF_Y,GF)
\qquad
(\phi:\sF_{HY}\Rightarrow F)\mapsto(G*\phi)\odot\sF_{\eta_Y}
$$
that assigns to every modification $\Theta:\phi\leadsto\phi'$
between pseudo-cones $\phi,\phi':\sF_{HY}\Rightarrow F$ the
modification $(G\circ\Theta)*\sF_{\eta_Y}$. We claim that $\mu$
is a quasi-inverse for $\theta^\dagger$. Indeed, for every
pseudo-cone $\psi:\sF_Y\Rightarrow GF$ we have
$\mu\circ\theta^\dagger(\psi)=
(G*((\eps*F)\odot(H*\psi)))\odot\sF_{\eta_Y}$, and example
\ref{ex_modifications}(i) and remark \ref{rem_construct-modifs}(iv)
yield an invertible modification
$$
\gamma^G_{\eps*F,H*\psi}*\sF_{\eta_Y}:
\mu\circ\theta^\dagger(\psi)\leadsto
(G*\eps*F)\odot(GH*\psi)\odot\sF_{\eta_Y}.
$$
Moreover, example \ref{ex_modifications}(ii) and remark
\ref{rem_construct-modifs}(iv) yield an invertible modification
$$
(G*\eps*F)*\tau^\eta_\psi:(G*\eps*F)\odot(GH*\psi)\odot\sF_{\eta_Y}
\leadsto(G*\eps*F)\odot(\eta*GF)\odot\psi
$$
which we may compose with the invertible modification
$$
(\Sigma\circ F)*\psi:
(G*\eps*F)\odot(\eta*GF)\odot\psi\leadsto\one_{GF}\odot\psi=\psi.
$$
Summing up, we obtain an invertible modification
$$
\Lambda_\psi:\mu\circ\theta^\dagger(\psi)\leadsto\psi
\qquad
i\mapsto(\Sigma_{Fi}*\psi_i)\odot(G(\eps_{Fi})*\tau^\eta_{\psi_i})
\odot(\gamma^G_{\eps_{Fi},H\psi_i}*\eta_Y)
$$
and we need to check that the rule $\psi\mapsto\Lambda_\psi$
yields a natural transformation
$$
\Lambda:\mu\circ\theta^\dagger\Rightarrow\one_{\sPsNat(\sF_Y,GF)}.
$$
Thus, let $\Xi:\psi\leadsto\psi'$ be any modification, and
notice that
$\mu\circ\theta^\dagger(\Xi)_i=G(\eps_{Fi}*H(\Xi_i))*\eta_Y$
for every $i\in\Ob(I)$; we compute :
$$
\begin{aligned}
\Xi_i\odot\Lambda_{\psi,i}=&\,
(\Sigma_{Fi}*\psi'_i)\odot(G(\eps_{Fi})*\eta_{GFi}*\Xi_i)\odot
(G(\eps_{Fi})*\tau^\eta_{\psi_i})\odot(\gamma^G_{\eps_{Fi},H\psi_i}*\eta_Y) \\
=&\,(\Sigma_{Fi}*\psi'_i)\odot(G(\eps_{Fi})*\tau^\eta_{\psi'_i})\odot
(G(\eps_{Fi})*GH(\Xi_i)*\eta_Y)\odot(\gamma^G_{\eps_{Fi},H\psi_i}*\eta_Y) \\
=&\,(\Sigma_{Fi}*\psi'_i)\odot(G(\eps_{Fi})*\tau^\eta_{\psi'_i})\odot
(\gamma^G_{\eps_{Fi},H\psi'_i}*\eta_Y)\odot(G(\eps_{Fi}*H(\Xi_i))*\eta_Y) \\
=&\,\Lambda_{\psi',i}\odot\mu\circ\theta^\dagger(\Xi)_i
\end{aligned}
$$
where the first identity follows from remark
\ref{rem_equiv-2-cat}(i), the second one from remark
\ref{rem_pseudo-natural}(i), and the third one from remark
\ref{rem_pseudo-funct}(ii). The assertion follows.

Lastly, we have $\theta^\dagger\circ\mu(\phi)=
(\eps*F)\odot H*((G*\phi)\odot\sF_{\eta_Y})$ for every
pseudo-cone $\phi:\sF_{HY}\Rightarrow F$, and we get an
invertible modification
$$
\Lambda'_\phi:\theta^\dagger\circ\mu(\phi)\leadsto
(\eps*F)\odot(HG*\phi)\odot\sF_{H\eta_Y}\leadsto
\phi\odot\sF_{\eps_{HY}}\odot\sF_{H\eta_Y}\leadsto\phi
$$
as the composition of $(\eps*F)*(\gamma^H_{\sF_{\eta_Y},G*\phi})^{-1}$
from example \ref{ex_modifications}(i), together with
$(\tau^\eps_\phi)^{-1}*\sF_{H\eta_Y}$ from example
\ref{ex_modifications}(ii) and $\phi*\sF_{\Sigma'_Y}$.
An explicit calculation as in the foregoing shows that
the rule $\phi\mapsto\Lambda'_\phi$ yields a natural
transformation $\Lambda:
\theta^\dagger\circ\mu\Rightarrow\one_{\sPsNat(\sF_{HY},F)}$
(details left to the reader), and concludes the proof
of (i).

Assertion (ii) follows from (i), by considering the
opposite $2$-categories : details left to the reader.
\end{proof}

\begin{proposition}\label{prop_2-cat-mitchell}
Let $\cA,\cB$ be two $2$-categories, $F:\cA\to\cB$ a fully
faithful pseudo-functor that admits a left $2$-adjoint
$G:\cB\to\cA$. The following holds :
\begin{enumerate}
\item
For every small $2$-category $I$ and every pseudo-functor
$\phi:I\to\cA$, if $F\circ\phi:I\to\cB$ admits a $2$-limit
$(L,\pi^L)$, then $\phi$ admits a $2$-limit $(M,\pi^M)$.
\item
In the situation of {\em (i)}, we have equivalences
$L\isom FM$ in $\cB$ and $M\isom GL$ in $\cA$.
\end{enumerate}
\end{proposition}
\begin{proof} Let $\eta:\one_\cB\Rightarrow FG$ and
$\eps:GF\Rightarrow\one_\cA$ be a unit and a counit for
the $2$-adjoint pair $(G,F)$, and set
$\pi^{GL}:=(\eps*\phi)\odot(G*\pi^L):\sF_{GL}\Rightarrow\phi$.
By the universality of $\pi^L$, there exist a $1$-cell
$f:FGL\to L$ in $\cB$ and an isomorphism of pseudo-cones
$\pi^L\odot\sF_f\isom F*\pi^{GL}$. Invoking example
\ref{ex_modifications}(i,ii) we deduce isomorphisms
of pseudo-cones :
$$
\begin{aligned}
\pi^L\!\odot\!\sF_{f\circ\eta_L}=
\pi^L\!\odot\!\sF_f\!\odot\!\sF_{\eta_L}\isom
(F*\pi^{GL})\!\odot\!(\eta*\sF_L)&\,\isom
(F*\eps*\phi)\!\odot\!(FG*\pi^L)\!\odot\!(\eta*\sF_L) \\
&\,\isom(F*\eps*\phi)\odot(\eta*F\phi)\circ\pi^L\isom\pi^L
\end{aligned}
$$
whence an isomorphism $f\circ\eta_L\isom\one_L$ in $\cB$, by
universality of $\pi^L$. On the other hand, since $F$ is fully
faithful, there exists a $1$-cell $h:GL\to GL$ in $\cA$ with
an isomorphism $Fh\isom\eta_L\circ f:FGL\to FGL$. Notice the
isomorphisms :
\set\begin{equation}\label{eq_like-mitchell}
Fh\circ\eta_L\isom\eta_L\circ f\circ\eta_L
\isom\eta_L\isom F\one_{GL}\circ\eta_L.
\end{equation}
On the other hand, using the pseudo-naturality of $\eps$ and
the triangular modifications for the pair $(\eta,\eps)$ we get
isomorphisms :
$$
\eps_{GL}\circ G(Fh\circ\eta_L)\isom\eps_{GL}\circ GFh\circ G\eta_L
\isom h\circ\eps_{GL}\circ G\eta_L\isom h
$$
and similarly :
$\eps_{GL}\circ G(F\one_{GL}\circ\eta_L)\isom\one_{GL}$. Combining with
\eqref{eq_like-mitchell}, we deduce an isomorphism $h\isom\one_{GL}$,
whence $\eta_L\circ f\isom F\one_{GL}\isom\one_{FGL}$, which shows that
$f$ is an equivalence in $\cB$. Especially, (ii) holds with $M:=GL$,
and it also follows that $(FGL,F*\pi^{GL})$ is a $2$-limit of
$F\circ\phi$. It remains to check that $(GL,\pi^{GL})$ is a $2$-limit
of $\phi$. To this aim, we consider for every $A\in\Ob(\cA)$ the
essentially commutative diagram of categories :
$$
\xymatrix{ \cA(A,GL) \ar[r] \ar[d]_{F_{A,GL}} &
\sPsNat(\sF_A,\phi) \ar[d]^U \\
\cB(FA,FGL) \ar[r] & \sPsNat(\sF_{FA},F\circ\phi)
}$$
whose top (resp. bottom) horizontal arrow is the pseudo-functor
induced by the pseudo-cone $\pi^{GL}$ (resp. $F*\pi^{GL}$) as in
definition \ref{def_pseudo-lim}(i). The pseudo-functor $U$ is
given by the rules : $\beta\mapsto F*\beta$ and
$\Theta\mapsto F*\Theta$ for every pseudo-cone
$\beta:\sF_A\Rightarrow\phi$ and every modification of pseudo-cones
$\Theta:\beta\leadsto\beta'$. Since $F$ is fully faithful, the
functor $F_{A,GL}$ is an equivalence, and the same holds for the
bottom horizontal arrow. We are then reduced to checking that
$U$ is an equivalence. To this aim, recall that $\eps$ is a
pseudo-natural equivalence (corollary \ref{cor_fully-faith-2-adjoint})
and fix a $1$-cell $\omega_A:A\to GFA$ in $\cA$ and an isomorphism
$$
\rho_A:\eps_A\circ\omega_A\isom\one_A.
$$
There follows an isomorphism $F\eps_A\circ F\omega_A\isom\one_{FA}$;
on the other hand, the triangular modifications for $(\eta,\eps)$
give an isomorphism $F\eps_A\circ\eta_{FA}\isom\one_{FA}$, so we get
as well an isomorphism
$$
\rho'_A:F\omega_A\isom\eta_{FA}.
$$
We consider the pseudo-functor
$$
V:\sPsNat(\sF_{FA},F\circ\phi)\to\sPsNat(\sF_A,\phi)
\qquad
\beta\mapsto(\eps*\phi)\odot(G*\beta)\odot\sF_{\omega_A}
$$
that assigns to every modification $\Theta:\beta\leadsto\beta'$
of pseudo-cones $\beta,\beta':\sF_{FA}\Rightarrow F\circ\phi$ the
modification $(\eps*\phi)*(G*\Theta)*\sF_{\omega_A}$. We construct
as follows isomorphisms of functors
$$
\Lambda:VU\isom\one_{\sPsNat(\sF_A,\phi)}
\qquad
\Lambda':UV\isom\one_{\sPsNat(\sF_{FA},F\circ\phi)}.
$$
For every pseudo-cone $\beta:\sF_A\Rightarrow\phi$ we have
$VU(\beta)=(\eps*\phi)\odot(GF*\beta)\odot\sF_{\omega_A}$. By
example \ref{ex_modifications}(ii), we have an invertible
modification
$$
\Xi_\beta:\beta\odot(\eps*\sF_A)\leadsto(\eps*\phi)\odot(GF*\beta)
\qquad
i\mapsto\tau^\eps_{\beta_i}
$$
where $\tau^\eps$ denotes the coherence constraint of $\eps$.
There follows an isomorphism $\Xi^{-1}_\beta*\sF_{\omega_A}:
VU(\beta)\isom\beta\odot(\eps*\sF_A)\odot\sF_{\omega_A}=
\beta\odot\sF_{\eps_A\circ\omega_A}$. We compose the latter with
the isomorphism
$\beta*\sF_{\rho_A}:\beta\odot\sF_{\eps_A\circ\omega_A}\isom\beta$
to obtain the isomorphism $\Lambda_\beta:VU(\beta)\isom\beta$.
In order to check the naturality of the rule :
$\beta\mapsto\Lambda_\beta$, we come down to verifying the
identity :
$$
((\eps*\phi)*(GF*\Theta))\odot\Xi_\beta=
\Xi_{\beta'}\odot(\Theta*(\eps*\sF_A))
$$
for every morphism of pseudo-cones $\Theta:\beta\leadsto\beta'$.
This follows from the naturality of $\tau^\eps$, applied to the
system of $2$-cells $\Theta_i:\beta_i\Rightarrow\beta'_i$ :
details left to the reader. Lastly, for every pseudo-cone
$\beta:\sF_{FA}\Rightarrow F\circ\phi$ we have
$UV(\beta)=F*((\eps*\phi)\odot(G*\beta)\odot\sF_{\omega_A})$.
By example \ref{ex_modifications}(i) we have the invertible
modification
$$
\Gamma_\beta:(F*\eps*\phi)\odot(FG*\beta)\odot\sF_{F\omega_A}
\leadsto F*((\eps*\phi)\odot(G*\beta)\odot\sF_{\omega_A})
\qquad
i\mapsto\gamma^F_{\omega_A,G\beta_i,\eps_{\phi(i)}}
$$
where $\gamma^F_{\bullet\bullet\bullet}$ is coherence constraint
of $F$ (for compositions of three $1$-cells). We compose
$\Gamma^{-1}_\beta$ with the invertible modification
$((F*\eps*\phi)\odot(FG*\beta))*\sF_{\rho'_A}$ to get
an isomorphism
$UV(\beta)\isom(F*\eps*\phi)\odot(FG*\beta)\odot(\eta*\sF_A)$.
Then, example \ref{ex_modifications}(ii) yields the invertible
modification
$$
\Xi'_\beta:(FG*\beta)\odot(\eta*\sF_A)\leadsto
(\eta*F*\phi)\odot\beta
\qquad
i\mapsto\tau^\eta_{\beta_i}
$$
where $\tau^\eta$ denotes the coherence constraint of $\eta$.
On the other hand we have the triangular modification
$\Sigma:(F*\eps)\odot(\eta*F)\leadsto\one_F$, so we can
compose further with
$$
((\Sigma*\phi)*\beta)\odot((F*\eps*\phi)*\Xi'_\beta):
(F*\eps*\phi)\odot(FG*\beta)\odot(\eta*\sF_A)\isom\beta
$$
to get the sought isomorphism $\Lambda'_\beta:UV(\beta)\isom\beta$.
Again, the naturality of the rule : $\beta\mapsto\Lambda'_\beta$
follows easily from the naturality of $\gamma^F_{\bullet\bullet\bullet}$
and $\tau^\eta$.
\end{proof}

\sset\subsubsection{}\label{subsec_2-final-pseudo-funct}
Let $I,J$ be two small $2$-categories, $\cC$ another $2$-category,
and $\phi:J\to I$, $F:I\to\cC$ two pseudo-functors. Let also
$(L,\pi)$ and $(L',\pi')$ be $2$-limits for $F$ and respectively
$F\circ\phi$. Then there exist a $1$-cell $f:L\to L'$ of $\cC$
and an invertible modification
$$
\Theta:\pi*\phi\leadsto\pi'\odot\sF_f.
$$
Moreover, if $f':L\to L'$ is another such $1$-cell with
an invertible modification
$\Theta':\pi*\phi\leadsto\pi'\odot\sF_{f'}$, there exists
a unique invertible $2$-cell in $\cC$ :
$$
\beta:f\isom f'
\qquad\text{such that}\qquad
\sF_\beta\odot\Theta=\Theta'.
$$
In particular, $f$ shall be an equivalence in $\cC$ if and only
if the same holds for $f'$ (remark \ref{rem_compose-equiv}(iii)),
so the property that $f$ is an equivalence is an intrinsic feature
of the pair $(\phi,F)$. Furthermore, $\Theta$ induces an essentially
commutative diagram of categories :
$$
{\diagram \cC(X,L) \ar[r]^-T \ar[d]_{f^X_*} &
\sPsNat(\sF_X,F) \ar[d]^{S^X} \\
\cC(X,L') \ar[r]^-{T'} & \sPsNat(\sF_X,F\circ\phi)
\enddiagram}
\qquad
\text{for every $X\in\Ob(\cC)$}
$$
whose horizontal arrows are the equivalences of definition
\ref{def_pseudo-lim}(i), where $f^X_*$ is defined as in example
\ref{ex_2-arrows}(iii), and $S^X$ is given by the rules :
$(\beta:\sF_X\Rightarrow F)\mapsto\beta*\phi$ and
$(\Xi:\beta\leadsto\beta')\mapsto\Xi\circ\phi$. Then the
required isomorphism of functors $S^X\circ T\isom T'\circ f^X_*$
is given by the rule : $(g:X\to L)\mapsto\Theta*\sF_g$. Hence,
$f$ is an equivalence in $\cC$ if and only if $S^X$ is an
equivalence of categories for every $X\in\Ob(\cC)$.

\begin{definition} We say that a pseudo-functor $\phi:J\to I$ is
{\em $2$-final} if for every $2$-complete $2$-category $\cC$ and
every pseudo-functor $F:I\to\cC$ the induced $1$-cell $f:L\to L'$
as in \eqref{subsec_2-final-pseudo-funct} is an equivalence. Likewise,
we say that $\phi$ is {\em $2$-cofinal} if $\phi^o:J^o\to I^o$ is
$2$-final.
\end{definition}

Unlike as in the case for usual category, we do not have a general
characterization of $2$-final or $2$-cofinal pseudo-functors, but
we can point out at least two simple, though useful, classes of
such pseudo-functors.
For our first example, let $\bone$ be a final object of $\bCat$,
{\em i.e.} a category with a unique object and a unique morphism
(then the unique morphism is the identity of the unique object);
let likewise $\bI$ be a final object of $2\tdu\bCat$, {\em i.e.}
a $2$-category with a unique object, a unique $1$-cell and a
unique $2$-cell. We say that a given object $i_0$ of $I$ is
{\em pseudo-initial} (resp. {\em pseudo-final}) if the category
$I(i_0,i)$ (resp. $I(i,i_0)$) is equivalent to $\bone$, for
every $i\in\Ob(I)$.

\begin{proposition}\label{prop_pseudo-initial}
Let $i_0$ be a pseudo-initial (resp. pseudo-final) object of
$I$, and $\phi:\bI\to I$ the strict pseudo-functor that maps
the object of\/ $\bI$ to $i_0$. Then $\phi$ is $2$-final
(resp. $2$-cofinal).
\end{proposition}
\begin{proof} Let $i_0$ be pseudo-initial; we need to check
that for every $2$-complete $2$-category $\cC$, every
pseudo-functor $F:I\to\cC$ and every $2$-limit $(L,\pi)$
of $F$, the pair $(L,\pi*\phi)$ is a $2$-limit of $F\circ\phi$,
and notice that the latter condition means that
$\pi_{i_0}:L\to L':=Fi_0$ is an equivalence. To this aim,
it suffices to construct a universal pseudo-cone
$\pi':\sF_{L'}\Rightarrow F$ with $\pi'_{i_0}=\one_{L'}$.
We may also assume that $F$ is unital, by lemma
\ref{lem_pseudo-trivial} and proposition
\ref{prop_towards-2-yoneda}.
Thus, let $i\in\Ob(I)$ be any object; by assumption,
there exists a $1$-cell $f_i:i_0\to i$ in $I$, and we
set $\pi'_i:=Ff_i:L'\to Fi$. Naturally, we take
$f_{i_0}:=\one_{i_0}$, so that $\pi'_{i_0}=\one_{L'}$, since
$F$ is unital. Next, for every one-cell $g:i\to i'$
in $I$ we have by assumption a unique $2$-cell
$\beta_g:g\circ f_i\Rightarrow f_{i'}$ in $I$, and it is
easily seen that $\beta_g$ must then be invertible (details
left to the reader); we set
$\tau_g:=(F\beta_g)\odot\gamma^F_{f_i,g}:
Ff\circ\pi_i\Rightarrow\pi'_{i'}$, where $\gamma^F$ denotes
the coherence constraint of $F$.

\begin{claim} The rule : $i\mapsto\pi'_i$ defines a pseudo-cone
$\pi':\sF_{L'}\Rightarrow F$ with coherence constraint given by
the system of $2$-cells $\tau_\bullet$.
\end{claim}
\begin{pfclaim} The uniqueness of $\beta_g$ implies that
$\beta_{\one_i}=\one_{f_i}$ for every $i\in\Ob(I)$ and
$$
\beta_h\odot(h*\beta_g)=\beta_{h\circ g}
\qquad
\text{for every pair of $1$-cells
$i\xrightarrow{g}i'\xrightarrow{h}i''$}.
$$
Since $F$ is unital, it follows already that
$\tau_{\one_i}=\one_{\pi'_i}$ for every $i\in\Ob(I)$. It remains
to check that
$$
\tau_h\odot(Fh*\tau_g)=\tau_{h\circ g}\odot(\gamma^F_{g,h}*\pi'_i)
$$
for every pair of $1$-cells $(g,h)$ as in the foregoing. We
compute :
$$
\begin{aligned}
\tau_h\odot(Fh*\tau_g)&\,=
F\beta_h\odot\gamma^F_{f_{i'},h}\odot(Fh*F\beta_g)
\odot(Fh*\gamma^F_{f_i,g}) \\
&\,=F\beta_h\odot F(h*\beta_g)\odot\gamma^F_{g\circ f_i,h}
\odot(Fh*\gamma^F_{f_i,g}) \\
&\,=F\beta_{h\circ g}\odot\gamma^F_{g\circ f_i,h}
\odot(Fh*\gamma^F_{f_i,g})
\end{aligned}
$$
whereas $\tau_{h\circ g}\odot(\gamma^F_{g,h}*\pi'_i)=F\beta_{g\circ h}
\odot\gamma^F_{f_i,g}\odot(\gamma^F_{g,h}*Ff_i)$, so the sought
identity follows from the composition axiom for $\gamma^F$.
\end{pfclaim}

To conclude, we need to check that functor
$$
\cC(X,L)\to\sPsNat(c_X,F)
\qquad
t\mapsto\pi'\odot\sF_t
$$
is an equivalence for every $X\in\Ob(\cC)$. The faithfulness
is clear, since $\pi'_{i_0}=\one_{L'}$. To check that the functor
is full, consider two $1$-cells $t,t':X\to L$ and a modification
$\Xi:\pi'\odot\sF_t\leadsto\pi'\odot\sF_{t'}$; we set
$\mu:=\Xi_{i_0}:t\Rightarrow t'$. Now, notice that
$\beta_{f_i}=\one_{f_i}$ for every $i\in\Ob(I)$, so that
$\tau_{f_i}=\one_{Ff_i}$, since $F$ is unital; then the compatibility
condition for $\Xi$, relative to the $1$-cell $f_i:i_0\to i$ yields:
$$
\Xi_i=Ff_i*\mu
\qquad
\text{for every $i\in\Ob(I)$}
$$
which shows that $\Xi=\pi'*\sF_\mu$, as required. Lastly, let
$\alpha:\sF_X\Rightarrow F$ be any pseudo-cone; set
$t:=\alpha_{i_0}:X\to L'$, and
$\Xi_i:=\tau^\alpha_{f_i}:Ff_i\circ t\Rightarrow\alpha_i$ for every
$i\in\Ob(I)$, where $\tau^\alpha$ is the coherence constraint of
$\alpha$. To conclude, it suffices to check that the rule :
$i\mapsto\Xi_i$ yields an invertible modification
$\Xi:\pi'\odot\sF_t\leadsto\alpha$. This comes down to showing :
$$
\tau^\alpha_g\odot(Fg*\tau^\alpha_{f_i})=
\tau^\alpha_{f_{i'}}\odot(\tau_g*\alpha_{i_0})
\qquad
\text{for every $1$-cell $g:i\to i'$ of $I$}.
$$
However, the left-hand side of the latter identity equals
$\tau^\alpha_{g\circ f_i}\odot(\gamma^F_{f_i,g}*\alpha_{i_0})$, by
the coherence axiom for $\tau^\alpha$, so we are reduced to
checking the identity :
$$
\tau^\alpha_{g\circ f_i}=\tau^\alpha_{f_{i'}}\odot(F\beta_g*\alpha_{i_0})
$$
which holds by naturality of $\tau^\alpha$.
\end{proof}

\sset\subsubsection{}\label{subsec_second-ex-of-2-cofinal}
For our second example, consider a small category $I$ (which
we regard as a $2$-category with trivial $2$-cells, as usual),
and a full subcategory $J$ of $I$, such that the inclusion
functor $\phi:J\to I$ admits a left adjoint $\sigma:I\to J$
with $\sigma\circ\phi=\one_J$.

\begin{proposition}\label{prop_another-2-cofinal}
In the situation of \eqref{subsec_second-ex-of-2-cofinal},
the pseudo-functor $\phi:J\to I$ is $2$-cofinal.
\end{proposition}
\begin{proof} For every pseudo-functor $F:I\to\cC$, every
pseudo-cocone $\pi:F\Rightarrow\sF_L$, and every $C\in\Ob(\cC)$
we have a commutative diagram of categories :
$$
\xymatrix{ & \cC(L,C) \ar[ld] \ar[rd] \\
\sPsNat(F,\sF_C) \ar[rr]^-{S_C} & & \sPsNat(F\circ\phi,\sF_C)
}$$
whose downward arrows are the functors $f\mapsto\pi\odot\sF_f$
and $f\mapsto(\pi*\phi)\odot\sF_f$ of definition
\ref{def_pseudo-lim}(ii), and where $S_C$ is given by the
rules : $(\beta:F\Rightarrow\sF_C)\mapsto\beta*\phi$ and
$(\Xi:\beta\leadsto\beta')\mapsto\Xi\circ\phi$. We come down
to checking that $S_C$ is an equivalence, for every such $F$.
To this aim, notice first that the adjoint pair $(\sigma,\phi)$
admits an adjunction whose counit $\sigma\circ\phi\Rightarrow\one_J$
is the identity natural transformation (proposition
\ref{prop_fullfaith-adjts}(iii)); in view of the triangular
identities of \eqref{subsec_adj-pair}, it follows that the
unit $\eta:\one_I\Rightarrow\phi\circ\sigma$ fulfills as well
the condition :
$$
\eta*\phi=\one_\phi.
$$
We consider for every $C\in\Ob(\cC)$ the functor
$$
S'_C:\sPsNat(F\circ\phi,\sF_C)\to\sPsNat(F,\sF_C)
\qquad
(\beta:F\circ\phi\Rightarrow\sF_C)\mapsto(\beta*\sigma)\odot(F*\eta)
$$
that assigns to every modification $\Xi:\beta\leadsto\beta'$
the modification $(\Xi\circ\sigma)*(F*\eta):
S'_C(\beta)\leadsto S'_C(\beta')$. Notice that for every
pseudo-natural transformation $\beta:F\circ\phi\Rightarrow\sF_C$
we have
$$
S_C\circ S'_C(\beta)=((\beta*\sigma)\odot(F*\eta))*\phi=
(\beta*\sigma*\phi)\odot(F*\eta*\phi)=\beta
$$
and likewise, if $\beta':F\circ\phi\Rightarrow\sF_C$ is another
such pseudo-natural transformation, and $\Xi:\beta\leadsto\beta'$
is any modification, we get $S_C\circ S'_C(\Xi)=\Xi$. On the other
hand, for every pseudo-natural transformation
$\beta:F\Rightarrow\sF_C$, example \ref{ex_modifications}(ii)
yields an invertible modification
$$
\tau^\beta_\eta:\beta=(\sF_C*\eta)\odot(\beta*\one_I)\leadsto
(\beta*\phi*\sigma)\odot(F*\eta)=S'_C\circ S_C(\beta).
$$
To conclude the proof, it suffices to check that the rule
$\beta\mapsto\tau^\beta_\eta$ yields a natural transformation
$\one_{\sPsNat(F,\sF_C)}\Rightarrow S'_C\circ S_C$. The latter
comes down to the commutativity of the diagram of modifications :
$$
\xymatrix{ \beta \ar@{~>}[r]^-{\tau^\beta_\eta} \ar@{~>}[d]_\Xi &
(\beta*\phi*\sigma)\odot(F*\eta)
\ar@{~>}[d]^{(\Xi\circ\phi\circ\sigma)*(F*\eta)} \\
\beta' \ar@{~>}[r]^-{\tau^{\beta'}_\eta} &
(\beta'*\phi*\sigma)\odot(F*\eta)
}$$
for every modification $\Xi:\beta\leadsto\beta'$. This in turn
amounts to the compatibility condition for $\Xi$, relative
to the $1$-cell $\eta_i:i\to\phi\circ\sigma(i)$, for every
$i\in\Ob(I)$.
\end{proof}

\subsection{\texorpdfstring{$2$}{2}-Categorical Kan extensions}
\label{sec_2-cat-Kan-ext}
Let $I,J,\cC$ be three $2$-categories, $\phi:I\to J$
a {\em unital} pseudo-functor, and suppose that $\cC$
is $2$-complete, $I$ is small and $J$ has small
$\Hom$-categories. We define as follows a pseudo-functor
$$
2\tdu\!\!\!\int_\phi:\sPsFun(I,\cC)\to\sPsFun(J,\cC).
$$
For every {\em unital} pseudo-functor $F:I\to\cC$ and every
$j\in\Ob(J)$, notice that the $2$-category $j/\phi I$ is
small under our assumptions, and choose a $2$-limit
$(L_{F,j},\pi^{F,j})$ of $F\circ\st_j:j/\phi I\to\cC$ (notation
of example \ref{ex_pseudo-functors}(ii)). Then, for every
$1$-cell $g:j\to j'$ in $J$ there exist a $1$-cell
$L_{F,g}:L_{F,j}\to L_{F,j'}$ in $\cC$ and an invertible modification
$$
\Omega^{F,g}:\pi^{F,j}*g^*\leadsto\pi^{F,j'}\odot\sF_{L_{F,g}}
$$
with the strict pseudo-functor $g^*:j'/\phi I\to j/\phi I$
defined as in example \ref{ex_pseudo-functors}(ii), and we set
\set\begin{equation}\label{eq_objs-and-1-cells}
2\tdu\!\!\!\int_\phi F(j):=L_{F,j}
\qquad
2\tdu\!\!\!\int_\phi F(g):=L_{F,g}.
\end{equation}
Obviously, for $g=\one_j$, we take $L_{F,g}$ to be the
identity $1$-cell of $L_{F,j}$, and $\Omega^{F,g}$ the
identity modification of $\pi^{F,j}$. If $h:j'\to j''$
is another $1$-cell in $J$, there exists a unique
invertible $2$-cell
$L_{F,g,h}:L_{F,h}\circ L_{F,g}\Rightarrow L_{F,h\circ g}$ in
$\cC$ such that
\set\begin{equation}\label{eq_triple-whammy}
(\pi^{F,j''}*\sF_{L_{F,g,h}})\odot(\Omega^{F,h}*\sF_{L_{F,g}})
\odot(\Omega^{F,g}\circ h^*)=\Omega^{F,h\circ g}.
\end{equation}
Moreover, if $g,g':j\to j'$ are two $1$-cells of $J$,
and $\beta:g\Rightarrow g'$ is a $2$-cell, we set
$$
\pi^{F,\beta}_f:=\tau^{F,j}_{\beta^*(f)}
\qquad
\text{for every $(f:j'\to\phi(i))\in\Ob(j'/\phi I)$}
$$
where $\tau^{F,j}$ denotes the coherence constraint of
$\pi^{F,j}$, and $\beta^*:g^*\Rightarrow g'^*$ is the strict
pseudo-natural transformation associated with $\beta$ as in
example \ref{ex_pseudo-functors}(ii). Since $F$ is unital,
$\pi^{F,\beta}_f$ is a $2$-cell
$\pi^{F,j}_{g^*(f)}\Rightarrow\pi^{F,j}_{g'^*(f)}$ in $\cC$, between
two $1$-cells $L_{F,j}\to Fi$. We claim that the rule :
$(f:j'\to\phi(i))\mapsto\pi^{F,\beta}_f$ defines an invertible
modification $\pi^{F,\beta}:\pi^{F,j}*g\leadsto\pi^{F,j}*g'$.
Indeed, notice first that for every pair of $1$-cells
$(i_1,g_1)\xrightarrow{(f,\alpha)}(i_2,g_2)\xrightarrow{(f',\alpha')}
(i_3,g_3)$ of $j/\phi(I)$ we have by definition :
\set\begin{equation}\label{eq_def-composition}
(f',\alpha')\circ(f,\alpha)=\alpha'\odot(\phi(f')*\alpha)
\odot((\gamma^\phi_{f,f'})^{-1}*g_1)
\qquad
\text{in $j/\phi(I)$}.
\end{equation}
Now, again using the unital assumption for $F$, the
assertion comes down to the identity :
$$
\tau^{F,j}_{g'^*(h,\alpha)}\odot(Fh*\tau^{F,j}_{\beta^*(f)})=
\tau^{F,j}_{\beta^*(f')}\odot\tau^{F,j}_{g^*(h,\alpha)}
$$
for every $1$-cell $(h,\alpha):(f:j'\to\phi(i))\to(f':j'\to\phi(i'))$
of $j'/\phi I$. Let us also remark that
$\beta^*(f')\circ g^*(\alpha,f)=g'^*(\alpha,f)\circ\beta^*(f)$,
and $\beta^*(f)=(\one_i,f*\beta)$ and likewise for $\beta^*(f')$;
since also $\phi$ is unital, we then have
$\gamma^\phi_{\one_i,h}=\one_h=\gamma^\phi_{h,\one_{i'}}$, and
the sought identity follows from \eqref{eq_def-composition}
and the coherence axiom for $\tau^{F,j}$. Consequently, there
exists a unique $2$-cell $L_{F,\beta}:L_{F,g}\Rightarrow L_{F,g'}$
in $\cC$ such that
\set\begin{equation}\label{eq_my-goodness}
\Omega^{F,g'}\odot\pi^{F,\beta}
=(\pi^{F,j'}*\sF_{L_{F,\beta}})\odot\Omega^{F,g}
\end{equation}
and we set
\set\begin{equation}\label{eq_define-for-2-cells}
2\tdu\!\!\!\int_\phi F(\beta):=L_{F,\beta}.
\end{equation}

\begin{lemma}\label{lem_construct-2-Fubini-F}
In the situation of \eqref{sec_2-cat-Kan-ext}, the system
of invertible $2$-cells $L_{F,\bullet\bullet}$ yields a coherence
constraint for a unital pseudo-functor $2\tdu\!\!\int_\phi F:J\to\cC$
defined on objects and $1$-cells by \eqref{eq_objs-and-1-cells},
and on $2$-cells by \eqref{eq_define-for-2-cells}.
\end{lemma}
\begin{proof} Let $\beta:g\Rightarrow g'$ and $\beta':g'\Rightarrow g''$
be two $2$-cells of $J$; a simple inspection shows that
$$
\beta'^*(f)\circ\beta^*(f)=(\beta'\odot\beta)^*(f)
\qquad
\text{for every $(f:j'\to\phi(i))\in\Ob(j'/\phi I)$}
$$
whence $\pi^{F,\beta'}\odot\pi^{F,\beta}=\pi^{F,\beta'\odot\beta}$.
On the other hand, we have
$$
\begin{aligned}
(\pi^{F,j'}*\sF_{L_{F,\beta'}})\odot(\pi^{F,j'}*\sF_{L_{F,\beta}})\odot
\Omega^{F,g}=&\,
(\pi^{F,j'}*\sF_{L_{F,\beta'}})\odot\Omega^{F,g'}\odot\pi^{F,\beta} \\
=&\,\Omega^{F,g''}\odot\pi^{F,\beta'}\odot\pi^{F,\beta}.
\end{aligned}
$$
Summing up, we already see that
$2\tdu\!\!\int_\phi F(\beta')\odot2\tdu\!\!\int_\phi F(\beta)=
2\tdu\!\!\int_\phi F(\beta'\odot\beta)$. Moreover, notice that
$\pi^{F,\one_g}$ is the identity automorphism of $\pi^{F,j}*g$
(remark \ref{rem_unital}(ii)), whence
$$
2\tdu\!\!\!\int_\phi F(\one_g)=\one_{L_{F,g}}
\qquad
\text{for every $1$-cell $g$ of $J$}.
$$
Next, let $g,g':j\to j'$ and $h,h':j'\to j''$ be four $1$-cells
of $J$ and $\beta:g\Rightarrow g'$, $\beta:h\Rightarrow h'$ two
$2$-cells; we need to check that
$$
\Bigl(2\tdu\!\!\!\int_\phi F(\beta'*\beta)\Bigr)\odot L_{F,g,h}=
L_{F,g',h'}\odot
\Bigl(2\tdu\!\!\!\int_\phi F(\beta')*2\tdu\!\!\!\int_\phi F(\beta)\Bigr)
$$
and in light of the foregoing, we may assume that either
$\beta=\one_g$ or $\beta'=\one_h$. In these cases, the assertion
comes down to the identities :
$$
\begin{aligned}
L_{F,h*\beta}\odot L_{F,g,h}=&\,L_{F,g',h}\odot(L_{F,h}*L_{F,\beta}) \\
L_{F,\beta'*g}\odot L_{F,g,h}=&\,L_{F,g,h'}\odot(L_{F,\beta'}*L_{F,g}).
\end{aligned}
$$
However, set
$X:=(\Omega^{F,h}*\sF_{L_{F,g}})\odot(\Omega^{F,g}\circ h^*)$. It
suffices to show that
$$
\begin{aligned}
Y:=(\pi^{F,j''}*(\sF_{L_{F,g',h}}\odot(\sF_{L_{F,h}}*\sF_{L_{F,\beta}}))\odot X
&\,=Z:=(\pi^{F,j''}*(\sF_{L_{F,h*\beta}}\odot\sF_{L_{F,g,h}}))\odot X \\
Y':=(\pi^{F,j''}*(\sF_{L_{F,g,h'}}\odot(\sF_{L_{F,\beta'}}*\sF_{L_{F,g}}))\odot X
&\,=Z':=(\pi^{F,j''}*(\sF_{L_{F,\beta'*g}}\odot\sF_{L_{F,g,h}}))\odot X.
\end{aligned}
$$
We compute :
$$
\begin{aligned}
Y=&\,(\pi^{F,j''}*\sF_{L_{F,g',h}})\odot(\Omega^{F,h}*\sF_{L_{F,g'}})\odot
((\pi^{F,j'}*h^*)*\sF_{L_{F,\beta}})\odot(\Omega^{F,g}\circ h^*) \\
=&\,(\pi^{F,j''}*\sF_{L_{F,g',h}})\odot(\Omega^{F,h}*\sF_{L_{F,g'}})\odot
(\Omega^{F,g'}\circ h^*)\odot(\pi^{F,\beta}\circ h^*) \\
=&\,\Omega^{F,h\circ g'}\odot(\pi^{F,\beta}\circ h^*).
\end{aligned}
$$
and on the other hand,
$Z=(\pi^{F,j''}*\sF_{L_{F,h*\beta}})\odot\Omega^{F,h\circ g}=
\Omega^{F,h\circ g'}\odot\pi^{F,\one_h*\beta}$. So, for the first
identity we are reduced to showing that
$\pi^{F,\beta}\circ h^*=\pi^{F,\one_h*\beta}$. The latter follows
by direct inspection. Next, we compute :
$$
\begin{aligned}
Y'=&\,(\pi^{F,j''}*\sF_{L_{F,g,h'}})\odot(\Omega^{F,h'}*\sF_{L_{F,g}})
\odot(\pi^{F,\beta'}*\sF_{L_{F,g}})\odot(\Omega^{F,g}\circ h^*) \\
=&\,(\Omega^{F,h'\circ g})\odot(\Omega^{F,g}*h'^*)^{-1}
\odot(\pi^{F,\beta'}*\sF_{L_{F,g}})\odot(\Omega^{F,g}\circ h^*)
\end{aligned}
$$
and on the other hand,
$Z'=(\pi^{F,j''}*\sF_{L_{F,\beta'*\one_g}})\odot\Omega^{F,h\circ g}=
\Omega^{F,h'\circ g}\odot\pi^{F,\beta'*g}$. So, for the second
identity we are reduced to checking that
$$
(\Omega^{F,g}*h'^*)\odot\pi^{F,\beta'*g}=
(\pi^{F,\beta'}*\sF_{L_{F,g}})\odot(\Omega^{F,g}\circ h^*).
$$
But this identity follows from the compatibility condition
for $\Omega^{F,g}$. Lastly, let us check the composition
and unit axioms for $2\tdu\!\!\int_\phi F$. Hence, let
$j\xrightarrow{g}j'\xrightarrow{h}j''\xrightarrow{k}j'''$
be three $1$-cells of $J$; we need to prove that
$$
L_{F,h\circ g,k}\odot(L_{F,k}*L_{F,g,h})=
L_{F,g,k\circ h}\odot(L_{F,h,k}*L_{F,g})
\qquad\text{and}\qquad
L_{F,h,\one_{j'}}=\one_{L_{F,h}}=L_{F,\one_j,h}.
$$
The two identities that express the unit axiom follow by a simple
inspection. For the composition axiom, set
$X:=(\Omega^{F,k}*\sF_{L_{F,h}}*\sF_{L_{F,g}})\odot
((\Omega^{F,h}\circ k^*)*\sF_{L_{F,g}})\odot(\Omega^{F,g}\circ h^*\circ k^*)$.
It suffices to show that
$$
Y:=(\pi^{F,j'''}*(\sF_{L_{F,h\circ g,k}}\odot(\sF_{L_{F,k}}*\sF_{L_{F,g,h}})))
\odot X=Z:=
(\pi^{F,j'''}*(\sF_{L_{F,g,k\circ h}}\odot(\sF_{L_{F,h,k}}*\sF_{L_{F,g}})))
\odot X.
$$
We compute :
$$
\begin{aligned}
Y&\!\!=\!(\pi^{F,j'''}\!\!\!*\!\sF_{L_{F,h\circ g,k}})\!\!\odot\!\!
(\Omega^{F,k}\!\!*\!\sF_{L_{F,h\circ g}})\!\!\odot\!\!
((\pi^{F,j''}\!\!*\!k^*)\!*\!\sF_{L_{F,g,h}})\!\!\odot\!\!
((\Omega^{F,h}\!\!\circ\!k^*)\!*\!\sF_{L_{F,g}})\!\!\odot\!\!
(\Omega^{F,g}\!\circ\!h^*\!\!\circ\!k^*) \\
&\!\!=\!(\pi^{F,j'''}*\sF_{L_{F,h\circ g,k}})\odot
(\Omega^{F,k}*\sF_{L_{F,h\circ g}})\odot(\Omega^{F,h\circ g}\circ k^*) \\
&\!\!=\!\Omega^{F,k\circ h\circ g}
\end{aligned}
$$
and a similar calculation yields as well $Z=\Omega^{F,k\circ h\circ g}$ :
details left to the reader.
\end{proof}

\begin{remark}\label{rem_strong-Fubini}
In the situation of \eqref{sec_2-cat-Kan-ext}, suppose
moreover that $\cC$ is strongly $2$-complete. In this case
we can choose for $(L_{F,j},\pi^{F,j})$ a strong $2$-limit of
$F\circ\st_j$, for every $j\in\Ob(J)$. Then, we may also choose
$L_{F,g}$ for every $1$-cell $g:j\to j'$ in $J$ such that
$\pi^{F,j}*g^*=\pi^{F,j'}\odot\sF_{L_{F,g}}$, and take $\Omega^{F,g}$
to be the identity modification. With these choices, it follows
easily that $L_{F,g,h}$ is an identity $2$-cell, for every composable
pair of $1$-cells $g,h$ of $J$. We conclude that, in this situation,
the pseudo-functor $2\tdu\!\!\int_\phi F$ of lemma
\ref{lem_construct-2-Fubini-F} shall be {\em strict}.
\end{remark}

\sset\subsubsection{}\label{subsec_2-Kan-pseudo-naturality}
Keep the notation of \eqref{sec_2-cat-Kan-ext}, and
let $F,F':I\to\cC$ be two unital pseudo-functors, and
$\beta:F\Rightarrow F'$ a pseudo-natural transformation.
Then, for every $j\in\Ob(J)$ there exist a $1$-cell
$L_{\beta,j}:L_{F,j}\to L_{F',j}$ and an invertible modification
$$
\Omega^{\beta,j}:
(\beta*\st_j)\odot\pi^{F,j}\leadsto\pi^{F',j}\odot\sF_{L_{\beta,j}}.
$$
Evidently, we take $L_{\one_F,j}:=\one_{L_{F,j}}$ and
$\Omega^{\one_F,j}:=\one_{\pi^{F,j}}$ for every such $F$ and $j$.
If $g:j\to j'$ is any $1$-cell of $J$, there exists a unique
invertible $2$-cell $L_{\beta,g}:L_{F',g}\circ L_{\beta,j}\Rightarrow
L_{\beta,j'}\circ L_{F,g}$ such that
\set\begin{equation}\label{eq_laughing-hard}
(\pi^{F',j'}*\sF_{L_{\beta,g}})\odot(\Omega^{F',g}*\sF_{L_{\beta,j}})\odot
(\Omega^{\beta,j}\circ g^*)=(\Omega^{\beta,j'}*\sF_{L_{F,g}})\odot
((\beta*\st_{j'})*\Omega^{F,g}).
\end{equation}
A direct inspection shows that $L_{\one_F,g}=\one_{L_{F,g}}$
for every such $F$ and $g$.

If $\beta':F'\Rightarrow F''$ is another pseudo-natural
transformation of unital pseudo-functors, there exists a
unique invertible $2$-cell
$L_{\beta,\beta',j}:L_{\beta',j}\circ L_{\beta,j}\Rightarrow
L_{\beta'\odot\beta,j}$ such that
\set\begin{equation}\label{eq_L_beta-beta-prime-j}
(\pi^{F'',j}*\sF_{L_{\beta,\beta',j}})\odot(\Omega^{\beta',j}*\sF_{L_{\beta,j}})
\odot((\beta'*\st_j)*\Omega^{\beta,j})=\Omega^{\beta'\odot\beta,j}.
\end{equation}

\begin{lemma}\label{lem_Fubini-on-pseudo-nats}
In the situation of \eqref{subsec_2-Kan-pseudo-naturality},
the system of $2$-cells $(L_{\beta,g}~|~g:j\to j')$ yields a
coherence constraint for a pseudo-natural transformation
$$
2\tdu\!\!\!\int_\phi\beta:
2\tdu\!\!\!\int_\phi F\Rightarrow 2\tdu\!\!\!\int_\phi F'
\qquad
j\mapsto L_{\beta,j}.
$$
\end{lemma}
\begin{proof} Clearly $L_{\beta,\one_j}=\one_{L_{\beta,j}}$ for every
$j\in\Ob(J)$. Next, let $j\xrightarrow{g}j'\xrightarrow{g'}j''$
be two $1$-cells of $J$; we need to check the identity :
$$
(L_{\beta,j''}*L_{F,g,g'})\odot(L_{\beta,g'}*L_{F,g})\odot
(L_{F',g'}*L_{\beta,g})=L_{\beta,g'\circ g}\odot(L_{F',g,g'}*L_{\beta,j}).
$$
To this aim, set
$$
\begin{aligned}
X&\,:=(\Omega^{F',g'}*\sF_{L_{F',g}}*\sF_{L_{\beta,j}})
\odot((\Omega^{F',g}\circ g'^*)*\sF_{L_{\beta,j}})\odot
(\Omega^{\beta,j}\circ g^*\circ g'^*) \\
Y&\,:=\pi^{F',j''}*((\sF_{L_{\beta,j''}}*\sF_{L_{F,g,g'}})\odot
(\sF_{L_{\beta,g'}}*\sF_{L_{F,g}})) \\
Z&\,:=Y\odot(\pi^{F',j''}*\sF_{L_{F',g'}}*\sF_{L_{\beta,g}})\odot X \\
U&\,:=(\pi^{F',j''}*(\sF_{L_{\beta,g'\circ g}}\odot
(\sF_{L_{F',g,g'}}*\sF_{L_{\beta,j}})))\odot X.
\end{aligned}
$$
It suffices to show that $U=Z$. We compute :
$$
\begin{aligned}
Z&=Y\!\odot\!(\Omega^{F',g'}\!\!*\!\sF_{L_{\beta,j'}}\!*\sF_{L_{F,g}})
\!\odot\!((\pi^{F',j'}\!\!*g'^*)*\sF_{L_{\beta,g}})\!\odot\!
((\Omega^{F',g}\circ g'^*)*\sF_{L_{\beta,j}})\!\odot\!
(\Omega^{\beta,j}\!\circ\!g^*\!\circ\!g'^*) \\
&=Y\odot(\Omega^{F',g'}*\sF_{L_{\beta,j'}}*\sF_{L_{F,g}})\odot
((\Omega^{\beta,j'}\circ g'^*)*\sF_{L_{F,g}})\odot
((\beta*\st_{j''})\!*\!(\Omega^{F,g}\circ g'^*)) \\
&=(\pi^{F',j''}\!\!*\!\sF_{L_{\beta,j''}*L_{F,g,g'}})
\!\!\odot\!\!(\Omega^{\beta,j''}\!\!*\!\sF_{L_{F,g'}*L_{F,g}})
\!\!\odot\!\!((\beta\!*\!\st_{j''})*\Omega^{F,g'}\!\!*\!\sF_{L_{F,g}})
\!\!\odot\!\!((\beta\!*\!\st_{j''})\!*\!(\Omega^{F,g}\!\circ\!g'^*)) \\
&=(\Omega^{\beta,j''}\!\!*\!\sF_{L_{F,g'\circ g}})
\!\!\odot\!\!((\beta\!*\!\st_{j''})\!*\!\pi^{F,j''}\!\!*\!\sF_{L_{F,g,g'}})
\!\!\odot\!\!((\beta\!*\!\st_{j''})\!*\!\Omega^{F,g'}\!\!*\!\sF_{L_{F,g}})
\!\!\odot\!\!((\beta\!*\!\st_{j''})\!*\!(\Omega^{F,g}\!\circ\!g'^*)) \\
&=(\Omega^{\beta,j''}*\sF_{L_{F,g'\circ g}})
\odot((\beta*\st_{j''})*\Omega^{F,g'\circ g})
\end{aligned}
$$
and on the other hand $U=(\pi^{F',j''}*\sF_{L_{\beta,g'\circ g}})
\odot(\Omega^{F',g'\circ g}*\sF_{L_{\beta,j}})\odot
(\Omega^{\beta,j}\circ g^*\circ g'^*)$, whence the contention.
Lastly, we check the naturality condition for $L_{\beta,g}$ : let
$f,g:j\to j'$ be two $1$-cells of $J$ and $\alpha:f\Rightarrow g$
a $2$-cell; we need to check the identity :
$$
L_{\beta,g}\odot(L_{F',\alpha}*L_{\beta,j})=
(L_{\beta,j'}*L_{F,\alpha})\odot L_{\beta,f}.
$$
To this aim, set $X:=(\Omega^{F',f}*\sF_{L_{\beta,j}})\odot
(\Omega^{\beta,j}\circ f^*)$; it suffices to show that
$$
Y':=(\pi^{F',j'}*(\sF_{L_{\beta,g}}\odot(\sF_{L_{F',\alpha}}*\sF_{L_{\beta,j}})))\odot X
=Z':=(\pi^{F',j'}*((\sF_{L_{\beta,j'}}*\sF_{L_{F,\alpha}})\odot\sF_{L_{\beta,f}}))\odot X.
$$

\begin{claim}\label{cl_isabeau}
$(\pi^{F',\alpha}*\sF_{L_{\beta,j}})\odot(\Omega^{\beta,j}\circ f^*)=
(\Omega^{\beta,j}\circ g^*)\odot((\beta*\st_{j'})*\pi^{F,\alpha})$.
\end{claim}
\begin{pfclaim} Recall that $\alpha$ induces a pseudo-natural
transformation $\alpha^*:f^*\Rightarrow g^*$, and
$\alpha^*(h)=(\one_i,\alpha):h\circ f\to h\circ g$ for every
$(h:j'\to\phi(i))\in\Ob(j'/\phi I)$. The compatibility condition
of $\Omega^{\beta,j}$ then yields a commutative diagram
$$
\xymatrix{ \beta_i\circ\pi^{F,j}_{f^*(h)}
\ar@{=>}[d]_{\beta_i*\pi^{F,\alpha}_h}
\ar@{=>}[rr]^-{\Omega^{\beta,j}_{f^*(h)}} & &
\pi^{F',j}_{f^*(h)}\circ L_{\beta,j} \ar@{=>}[d]^{\pi^{F',\alpha}_h} \\
\beta_i\circ\pi^{F,j}_{g^*(h)}
\ar@{=>}[rr]^-{\Omega^{\beta,j}_{g^*(h)}} & &
\pi^{F',j}_{g^*(h)}\circ L_{\beta,j}.
}$$
The claim is an immediate consequence.
\end{pfclaim}

Using claim \ref{cl_isabeau}, we may compute :
$$
\begin{aligned}
Y'&\,=(\pi^{F',j'}*\sF_{L_{\beta,g}})
\odot(\Omega^{F',g}*\sF_{L_{\beta,j}})\odot(\pi^{F',\alpha}*\sF_{L_{\beta,j}})
\odot(\Omega^{\beta,j}\circ f^*) \\
&\,=(\pi^{F',j'}*\sF_{L_{\beta,g}})
\odot(\Omega^{F',g}*\sF_{L_{\beta,j}})\odot(\Omega^{\beta,j}\circ g^*)
\odot((\beta*\st_{j'})*\pi^{F,\alpha}) \\
&\,=(\Omega^{\beta,j'}*\sF_{L_{F,g}})\odot((\beta*\st_{j'})*\Omega^{F,g})
\odot((\beta*\st_{j'})*\pi^{F,\alpha}) \\
&\,=(\Omega^{\beta,j'}*\sF_{L_{F,g}})\odot
(\beta*\st_{j'}*((\pi^{F,j'}*\sF_{L_{F,\alpha}})\odot\Omega^{F,f}))
\end{aligned}
$$
and on the other hand $Z'=(\pi^{F',j'}*(\sF_{L_{\beta,j'}}*\sF_{L_{F,\alpha}}))
\odot(\Omega^{\beta,j'}*\sF_{L_{F,f}})\odot((\beta*\st_{j'})*\Omega^{F,f})$.
So, we are reduced to proving that
$$
(\Omega^{\beta,j'}*\sF_{L_{F,g}})\odot
(\beta*\st_{j'}*\pi^{F,j'}*\sF_{L_{F,\alpha}})=
(\pi^{F',j'}*(\sF_{L_{\beta,j'}}*\sF_{L_{F,\alpha}}))
\odot(\Omega^{\beta,j'}*\sF_{L_{F,f}}).
$$
But the latter follows as usual from remark \ref{rem_equiv-2-cat}(i).
\end{proof}

\begin{lemma}\label{lem_coherence-for-Fubini}
In the situation of \eqref{subsec_2-Kan-pseudo-naturality}, the
system of $2$-cells $(L_{\beta,\beta',j}~|~j\in\Ob(J))$ yields an
invertible modification
$$
L_{\beta,\beta'}:2\tdu\!\!\!\int_\phi\beta'\odot
2\tdu\!\!\!\int_\phi\beta\leadsto 2\tdu\!\!\!\int_\phi\beta'\odot\beta.
$$
\end{lemma}
\begin{proof} Let $g:j\to j'$ be any $1$-cell of $J$; we need
to show the identity
$$
(L_{\beta,\beta',j'}*L_{F,g})\odot(L_{\beta',j'}*L_{\beta,g})\odot
(L_{\beta',g}*L_{\beta,j})=L_{\beta'\odot\beta,g}\odot(L_{F'',g}*L_{\beta,\beta',j}).
$$
To this aim, let us set
$X:=((\beta'*\st_{j'})*\Omega^{\beta,j'}*\sF_{L_{F,g}})^{-1}
\odot(\Omega^{\beta',j'}*\sF_{L_{\beta,j'}}*\sF_{L_{F,g}})^{-1}$ and
$$
\begin{aligned}
Y_1&\,:=X\odot(\pi^{F'',j'}*\sF_{L_{\beta',j'}}*\sF_{L_{\beta,g}}) \\
Y_2&\,:=\pi^{F'',j'}*((\sF_{L_{\beta',g}}*\sF_{L_{\beta,j}})\odot
(\sF_{L_{F'',g}}*\sF_{L_{\beta,\beta',j}})^{-1}) \\
Z&\,:=X\odot(\pi^{F'',j'}*((\sF_{L_{\beta,\beta',j'}}*\sF_{L_{F,g}})^{-1}
\odot\sF_{L_{\beta'\odot\beta,g}})).
\end{aligned}
$$
It suffices to check that $Y:=Y_1\odot Y_2=Z$.
Set as well $U:=(((\beta'\odot\beta)*\st_{j'})*\Omega^{F,g})$;
we compute :
$$
\begin{aligned}
Z&\,=(\Omega^{\beta'\odot\beta,j'}*\sF_{L_{F,g}})^{-1}\odot
(\pi^{F'',j'}*\sF_{L_{\beta'\odot\beta,g}}) \\
&\,=U\odot(\Omega^{\beta'\odot\beta,j}\circ g^*)^{-1}\odot
(\Omega^{F'',g}*\sF_{L_{\beta'\odot\beta,j}})^{-1}.
\end{aligned}
$$
On the other hand, set $U':=(\Omega^{F'',g}*\sF_{L_{\beta',j}})^{-1}
\odot(\pi^{F'',j'}*\sF_{L_{\beta',g}})^{-1})*\sF_{L_{\beta,j}}$; then
$$
\begin{aligned}
Y&\,=((\beta'*\st_{j'})*\Omega^{\beta,j'}*\sF_{L_{F,g}})^{-1}
\!\odot\!((\beta'*\st_{j'})*\pi^{F',j'}*\sF_{L_{\beta,g}})
\!\odot\!(\Omega^{\beta',j'}*\sF_{L_{F',g}}*\sF_{L_{\beta,j}})^{-1}
\!\odot\!Y_2 \\
&\,=U\odot((\beta'*\st_{j'})*((\Omega^{\beta,j}\circ g^*)^{-1}
\odot(\Omega^{F',g}*\sF_{L_{\beta,j}})^{-1})\odot
(\Omega^{\beta',j'}*\sF_{L_{F',g}}*\sF_{L_{\beta,j}})^{-1}\odot Y_2 \\
&\,=U\odot((\beta'*\st_{j'})*(\Omega^{\beta,j}\circ g^*)^{-1})
\odot((\Omega^{\beta',j}\circ g^*)^{-1}*\sF_{L_{\beta,j}})
\odot U'\odot Y_2.
\end{aligned}
$$
Therefore, we are reduced to checking the identity :
$$
((\pi^{F'',j}*g^*)*\sF_{\beta,\beta',j})^{-1}\odot
(\Omega^{F'',g}*\sF_{L_{\beta'\odot\beta,j}})^{-1}=U'\odot Y_2.
$$
However, we have :
$$
((\pi^{F'',j}*g^*)*\sF_{\beta,\beta',j})^{-1}\odot
(\Omega^{F'',g}*\sF_{L_{\beta'\odot\beta,j}})^{-1}\!=\!
(\Omega^{F'',g}*\sF_{L_{\beta',j}*L_{\beta,j}})^{-1}
\odot(\pi^{F'',j'}*\sF_{L_{F'',g}}*\sF_{\beta,\beta',j})^{-1}
$$
so we are further reduced to showing that :
$$
(\pi^{F'',j'}*\sF_{L_{F'',g}}*\sF_{\beta,\beta',j})^{-1}=
(\pi^{F'',j'}*\sF_{L_{\beta',g}}*\sF_{L_{\beta,j}})^{-1}\odot Y_2
$$
which is obvious.
\end{proof}

\sset\subsubsection{}\label{subsec_Fubini-mods}
Let $F,F':I\to\cC$ be two unital pseudo-functors,
$\beta,\beta':F\Rightarrow F'$ two pseudo-natural transformations,
and $\Xi:\beta\leadsto\beta'$ a modification. Then for every
$j\in\Ob(J)$ there exists a unique modification
$L_{\Xi,j}:L_{\beta,j}\leadsto L_{\beta',j}$ such that
$$
\Omega^{\beta',j}\odot((\Xi\circ\st_j)*\pi^{F,j})=
(\pi^{F',j}*\sF_{L_{\Xi,j}})\odot\Omega^{\beta,j}.
$$

\begin{lemma} In the situation of \eqref{subsec_Fubini-mods},
the rule $j\mapsto L_{\Xi,j}$ defines a modification
$$
2\tdu\!\!\!\int_\phi\Xi:2\tdu\!\!\!\int_\phi\beta\leadsto
2\tdu\!\!\!\int_\phi\beta'.
$$
\end{lemma}
\begin{proof} Let $g:j\to j'$ be any $1$-cell of $J$; we need
to check the identity :
$$
(L_{\Xi,j'}*L_{F,g})\odot L_{\beta,g}=L_{\beta',g}\odot(L_{F',g}*L_{\Xi,j})
$$
and as usual, it suffices to show that
$$
Y:=\pi^{F',j'}*((\sF_{L_{\Xi,j'}}*\sF_{L_{F,g}})\odot\sF_{L_{\beta,g}})=
Z:=\pi^{F',j'}*(\sF_{L_{\beta',g}}\odot(\sF_{L_{F',g}}*\sF_{L_{\Xi,j}})).
$$
We compute :
$$
\begin{aligned}
Y&\!=\!((\Omega^{\beta',j'}\odot((\Xi\circ\st_{j'})*\pi^{F,j'})\odot
(\Omega^{\beta,j'})^{-1})*\sF_{L_{F,g}})\odot(\pi^{F',j'}*\sF_{L_{\beta,g}}) \\
&\!=\!((\Omega^{\beta',j'}\!\odot\!
((\Xi\circ\st_{j'})*\pi^{F,j'}))*\sF_{L_{F,g}})\!\odot\!
((\beta\!*\!\st_{j'})*\Omega^{F,g})\!\odot\!(\Omega^{\beta,j}\circ g^*)^{-1}
\!\odot\!(\Omega^{F',g}*\sF_{L_{\beta,j}})^{-1} \\
&\!=\!(\Omega^{\beta',j'}\!\!*\!\sF_{L_{F,g}})\!\odot\!
((\beta'*\st_{j'})\!*\!\Omega^{F,g})
\!\odot\!((\Xi\circ\st_{j'})*(\pi^{F,j}\!*\!g^*))\!\odot\!
(\Omega^{\beta,j}\!\circ\!g^*)^{-1}\!\odot\!
(\Omega^{F',g}\!*\!\sF_{L_{\beta,j}})^{-1} \\
&\!=\!(\Omega^{\beta',j'}\!\!*\!\sF_{L_{F,g}})\!\odot\!
((\beta'*\st_{j'})\!*\!\Omega^{F,g})
\!\odot\!(\Omega^{\beta',j}\!\circ\!g^*)^{-1}\!\odot\!
((\pi^{F',j}*g^*)*\sF_{L_{\Xi,j}})\!\odot\!
(\Omega^{F',g}*\sF_{L_{\beta,j}})^{-1} \\
&\!=\!(\pi^{F',j'}*\sF_{L_{\beta',g}})\odot(\Omega^{F',g}*\sF_{L_{\beta',j}})
\odot((\pi^{F',j}*g^*)*\sF_{L_{\Xi,j}})\!\odot\!
(\Omega^{F',g}*\sF_{L_{\beta,j}})^{-1} \\
&\!=\!(\pi^{F',j'}*\sF_{L_{\beta',g}})\odot(\Omega^{F',g}*\sF_{L_{\beta',j}})
\odot(\Omega^{F',g}*\sF_{\beta',j})^{-1}\odot
(\pi^{F',j'}*\sF_{L_{F,g}}*\sF_{\Xi,j}) \\
&\!=\!(\pi^{F',j'}*\sF_{L_{\beta',g}})\odot(\pi^{F',j'}*\sF_{L_{F,g}}*\sF_{\Xi,j})
\end{aligned}
$$
whence the contention.
\end{proof}

\begin{proposition}\label{prop_2-Fubini}
Let $I,J,\cC$ be as in \eqref{sec_2-cat-Kan-ext},
and $\phi:I\to J$ any pseudo-functor. The rules:
$$
F\mapsto 2\tdu\!\!\!\int_{\phi^u} F^u
\qquad
\beta\mapsto 2\tdu\!\!\!\int_{\phi^u}\beta^u
\qquad
\Xi\mapsto 2\tdu\!\!\!\int_{\phi^u}\Xi^u
$$
for every pseudo-functor $F:I\to\cC$, every pseudo-natural
transformation $\beta$ between pseudo-functors $F,F':I\to\cC$,
and every modification $\Xi$ of such pseudo-natural transformations
(notation of remark {\em\ref{rem_uniPsFun}(i)}) define a unital
pseudo-functor
$$
2\tdu\!\!\!\int_\phi:\sPsFun(I,\cC)\to\sPsFun(J,\cC)
$$
with coherence constraint given by the system of invertible
modifications $L_{\beta,\beta'}$ of lemma
{\em\ref{lem_coherence-for-Fubini}}. We call this pseudo-functor
the {\em right $2$-Kan extension along $\phi$}.
\end{proposition}
\begin{proof} The sought functor shall be the composition
of the strict $2$-equivalence of remark \ref{rem_uniPsFun}(i)
and a similar pseudo-functor
$\uniPsFun(I,\cC)\to\uniPsFun(J,\cC)$ given by the foregoing
rules; we are therefore reduced to prove the existence of the
latter. Now, a direct inspection shows that
$L_{\one_F,\beta,j}=\one_{L_{\beta,j}}=L_{\beta,\one_{F'},j}$ for every
pair of unital pseudo-functors $F,F':I\to\cC$, every
$\beta:F\Rightarrow F'$ and every $j\in\Ob(J)$, whence the
unit axiom for the system of modifications $L_{\beta,\beta'}$.
In order to check the composition axiom, we need to show
the identity :
$$
L_{\beta'\odot\beta,\beta'',j}\odot(L_{\beta'',j}*L_{\beta,\beta',j})=
L_{\beta,\beta''\odot\beta',j}\odot(L_{\beta',\beta'',j}*L_{\beta,j})
$$
for every $j\in\Ob(J)$ and every three pseudo-natural
transformations of unital pseudo-functors
$\beta:F\Rightarrow F'$, $\beta':F'\Rightarrow F''$,
$\beta'':F''\Rightarrow F'''$. To this aim, set
$X:=\Omega^{\beta''\odot\beta'\odot\beta,j}$, as well as :
$$
\begin{aligned}
Y&\,:=(\pi^{F''',j}*((\sF_{L_{\beta'',j}}*\sF_{L_{\beta,\beta',j}})^{-1}\odot
\sF_{L_{\beta'\odot\beta,\beta'',j}}^{-1}))\odot X \\
Z&\,:=(\pi^{F''',j}*((\sF_{L_{\beta',\beta'',j}}*\sF_{L_{\beta,j}})^{-1}\odot
\sF_{L_{\beta,\beta''\odot\beta',j}}^{-1}))\odot X.
\end{aligned}
$$
It suffices to prove that $Y=Z$. We compute :
$$
\begin{aligned}
Y&\,=(\pi^{F''',j}*\sF_{L_{\beta'',j}}*\sF_{L_{\beta,\beta',j}})^{-1}\odot
(\Omega^{\beta'',j}*\sF_{L_{\beta'\odot\beta,j}})\odot
((\beta''*\st_j)*\Omega^{\beta'\odot\beta,j}) \\
&\,=(\Omega^{\beta'',j}*\sF_{L_{\beta',j}}*\sF_{L_{\beta,j}})\odot
((\beta''*\st_j)*\pi^{F'',j}*\sF_{L_{\beta,\beta',j}})^{-1}\odot
((\beta''*\st_j)*\Omega^{\beta'\odot\beta,j}) \\
&\,=(\Omega^{\beta'',j}*\sF_{L_{\beta',j}}*\sF_{L_{\beta,j}})\odot
((\beta''*\st_j)*((\Omega^{\beta',j}*\sF_{L_{\beta,j}})
\odot((\beta'*\st_j)*\Omega^{\beta,j})))
\end{aligned}
$$
and likewise :
$$
\begin{aligned}
Z&\,=(\pi^{F''',j}*\sF_{L_{\beta',\beta'',j}}*\sF_{L_{\beta,j}})^{-1}\odot
(\Omega^{\beta''\odot\beta',j}*\sF_{L_{\beta,j}})
\odot((\beta''\odot\beta'*\st_j)*\Omega^{\beta,j}) \\
&\,=(((\Omega^{\beta'',j}*\sF_{L_{\beta',j}})
\odot((\beta''*\st_j)*\Omega^{\beta',j}))*\sF_{L_{\beta,j}})\odot
((\beta''\odot\beta'*\st_j)*\Omega^{\beta,j})
\end{aligned}
$$
whence the contention. Next, consider three pseudo-natural
transformations $\beta,\beta',\beta'':F\Rightarrow F'$ and
two modifications $\Xi:\beta\leadsto\beta'$ and
$\Xi':\beta'\leadsto\beta''$. For every $j\in\Ob(J)$ we have :
$$
\Omega^{\beta'',j}\odot(((\Xi'\odot\Xi)\circ\st_j)*\pi^{F,j})\!=\!
(\pi^{F',j}*\sF_{L_{\Xi',j}})\odot\Omega^{\beta',j}\odot
((\Xi\circ\st_j)*\pi^{F,j})\!=\!(\pi^{F',j}*\sF_{L_{\Xi',j}}*\sF_{L_{\Xi,j}})
\odot\Omega^{\beta,j}
$$
whence $L_{\Xi'\odot\Xi,j}=L_{\Xi',j}\odot L_{\Xi,j}$. This shows that
$2\tdu\!\!\int_\phi(\Xi'\odot\Xi)=
2\tdu\!\!\int_\phi\Xi'\odot2\tdu\!\!\int_\phi\Xi$. Moreover, a
simple inspection shows that $L_{\one_\beta,j}=\one_{L_{\beta,j}}$
for every pseudo-natural transformation $\beta$ and every
$j\in\Ob(J)$, so $2\tdu\!\!\int_\phi(\one_\beta)$ is the
identity modification of $2\tdu\!\!\int_\phi\beta$.

Lastly, let $F,F',F'':I\to\cC$ be three unital pseudo-functors,
$\beta,\beta:F\Rightarrow F'$ and $\alpha,\alpha':F'\Rightarrow F''$
four pseudo-natural transformations, and $\Xi:\beta\leadsto\beta'$,
$\Theta:\alpha\leadsto\alpha'$ two modifications; we need to show that
$$
L_{\Theta*\Xi,j}\odot L_{\beta,\alpha,j}=
L_{\beta',\alpha',j}\odot(L_{\Theta,j}*L_{\Xi,j})
\qquad
\text{for every $j\in\Ob(J)$}
$$
and by the foregoing, we may assume that either
$\Xi=\one_\beta$ or $\Theta=\one_\alpha$. Suppose first that
$\Theta=\one_\alpha$, and let $X:=(\Omega^{\alpha,j}*\sF_{\beta,j})
\odot((\alpha*\st_j)*\Omega^{\beta,j})$; it suffices to prove that
$$
Y:=(\pi^{F'',j}*(\sF_{L_{\alpha*\Xi,j}}\odot\sF_{L_{\beta,\alpha,j}}))\odot X=
Z:=(\pi^{F'',j}*(\sF_{L_{\beta',\alpha,j}}\odot
(\sF_{L_{\alpha,j}}*\sF_{L_{\Xi,j}})))\odot X.
$$
We compute :
$$
\begin{aligned}
Z&\,=(\pi^{F'',j}*\sF_{L_{\beta',\alpha,j}})\odot(\Omega^{\alpha,j}*\sF_{\beta',j})
\odot((\alpha*\st_j)*\pi^{F',j}*\sF_{L_{\Xi,j}})\odot
((\alpha*\st_j)*\Omega^{\beta,j}) \\
&\,=(\pi^{F'',j}*\sF_{L_{\beta',\alpha,j}})\odot(\Omega^{\alpha,j}*\sF_{\beta',j})
\odot((\alpha*\st_j)*(\Omega^{\beta',j}\odot((\Xi\circ\st_j)*\pi^{F,j}))) \\
&\,=\Omega^{\alpha\odot\beta',j}\odot
((\alpha*\st_j)*((\Xi\circ\st_j)*\pi^{F,j})))
\end{aligned}
$$
and on the other hand,
$Y=(\pi^{F'',j}*\sF_{L_{\alpha*\Xi,j}})\odot\Omega^{\alpha\odot\beta}$, whence
the contention. A similar calculation, left to the reader, settles
the case where $\Xi=\one_\beta$, and concludes the proof.
\end{proof}

\sset\subsubsection{}\label{subsec_left-2-Kan-ext}
On the other hand, if $\cC$ is $2$-cocomplete, then $\cC^o$ is
$2$-complete (remark \ref{rem_pseudo-limit}(vi)), and in view
of the strict isomorphisms of \eqref{subsec_opposite-mods},
we get a {\em left $2$-Kan extension along $\phi$} :
$$
2\tdu\!\!\!\int^\phi:\sPsFun(I,\cC)\to\sPsFun(J,\cC)
\qquad
F\mapsto\Bigl(2\tdu\!\!\!\int_{\phi^o}F^o\Bigr)^o.
$$
Explicitly, for every $j\in\Ob(J)$, the object
$(2\tdu\!\!\int^\phi F)(j)$ represents the $2$-colimit of
the pseudo-functor
$F\circ{}^o\ss_{\,{}^o\!j}:{}^o({}^o\phi({}^oI)/\,{}^o\!j)\to\cC$.

\begin{remark}\label{rem_strong-2-Kan-ext}
(i)\ \
In the situation of \eqref{subsec_2-Kan-pseudo-naturality},
suppose that $\cC$ is strongly $2$-complete, and choose strong
$2$-limits $(L_{F,j},\pi^{F,j})$, identity modifications $\Omega^{F,g}$,
and $1$-cells $L_{F,g}$ as in remark \ref{rem_strong-Fubini}.
Then we may further choose $L_{\beta,j}$ such that
$(\beta*\st_j)\odot\pi^{F,j}=\pi^{F',j}\odot\sF_{L_{\beta,j}}$, and
let $\Omega^{\beta,j}$ as well be the corresponding identity
modification. With such choices, it follows easily that
$L_{\beta,g}$ shall be an identity $2$-cell, for every $1$-cell
$g$ of $J$, and the same for $L_{\beta,\beta'}$. Therefore, in
this situation the pseudo-natural transformation
$2\tdu\!\!\int_\phi\beta$ of lemma \ref{lem_Fubini-on-pseudo-nats}
shall be {\em strict} as well, and moreover
$(2\tdu\!\!\int_\phi\beta')\odot(2\tdu\!\!\int_\phi\beta)=
2\tdu\!\!\int_\phi(\beta'\odot\beta)$.

(ii)\ \
We conclude that, with these choices, the pseudo-functor
$2\tdu\!\!\int_\phi$ of proposition \ref{prop_2-Fubini}
factors through a {\em strict} pseudo-functor called the
{\em strong right $2$-Kan extension along $\phi$} :
$$
\sPsFun(I,\cC)\to\stPsFun(J,\cC)
$$
(notation of definition \ref{def_PsNat}(iii)).
Likewise, if $\cC$ is strongly $2$-cocomplete, we obtain
as in \eqref{subsec_left-2-Kan-ext}, a {\em strong left
$2$-Kan extension along $\phi$} which is another strict
pseudo-functors with values in $\stPsFun(J,\cC)$, namely
the opposite of the strong right $2$-Kan extension along
$\phi^o$.
\end{remark}

\sset\subsubsection{}\label{subsec_build-adjoint}
Keep the notation of \eqref{sec_2-cat-Kan-ext}, and suppose
now that {\em all the $2$-cells of $J$ are invertible}. Under
this assumption, we easily see that the rule :
$(i,f:j\to\phi(i))\mapsto f$ yields a pseudo-cone
$$
\phi^*_j:\sF_j\Rightarrow\phi\circ\st_j
\qquad
\text{for every $j\in\Ob(J)$}
$$
whose coherence constraint is given by the rule :
$(h,\alpha)\mapsto\alpha$ (notation of example
\ref{ex_pseudo-functors}(ii); notice that this would not be
a coherence constraint, if $\alpha$ were not invertible :
details left to the reader). Moreover, we have :
\set\begin{equation}\label{eq_box-or-place}
\phi^*_{j'}\odot\sF_g=\phi^*_j*g^*
\qquad
\text{for every $1$-cell $g:j\to j'$ of $J$}.
\end{equation}
Furthermore, for every pair of $1$-cells $g,g':j\to j'$ and
every $2$-cell $\alpha:g\Rightarrow g'$ we get a modification
$$
\phi^*_j*\alpha^*:\phi^*_j*g^*\leadsto\phi^*_j*g'^*
\qquad
(f:j'\to\phi(i))\mapsto(\alpha*f:g^*(f)\Rightarrow g'^*(f)).
$$
Let $F:I\to\cC$ and $G:J\to\cC$ be two {\em unital}
pseudo-functors, and $\beta:G\circ\phi\Rightarrow F$ a
pseudo-natural transformation; there exist a $1$-cell
$\beta^\dagger_j:Gj\to L_{F,j}$ and an invertible modification
$$
\Theta^{\beta,j}:(\beta*\st_j)\odot(G*\phi^*_j)\leadsto
\pi^{F,j}\odot\sF_{\beta^\dagger_j}
\qquad
\text{for every $j\in\Ob(J)$}.
$$
Next, let $g:j\to j'$ be any $1$-cell of $J$, and notice
that $G*\phi^*_j*g^*=G*(\phi^*_{j'}\odot\sF_g)$, due to
\eqref{eq_box-or-place}; according to example
\ref{ex_modifications}(i) the coherence constraint $\gamma^G$
of $G$ induces an invertible modification
$$
\Gamma^{G,g}:(G*\phi^*_{j'})\odot\sF_{Gg}\leadsto G*\phi^*_j*g^*
\qquad
(f:j'\to\phi(i))\mapsto\gamma^G_{g,f}.
$$
Then there exists a unique invertible $2$-cell
$\tau^{\beta^\dagger}_g:L_{F,g}\circ\beta^\dagger_j
\Rightarrow\beta^\dagger_{j'}\circ Gg$
such that :
\set\begin{equation}\label{eq_redo-sat-night}
(\pi^{F,j'}*\sF_{\tau^{\beta^\dagger}_g})\odot(\Omega^{F,g}*\sF_{\beta^\dagger_j})
\odot(\Theta^{\beta,j}\circ g^*)\odot((\beta*\st_{j'})*\Gamma^{G,g})
=\Theta^{\beta,j'}*\sF_{Gg}.
\end{equation}

\begin{lemma}\label{lem_build-2-adj-step-1}
With the notation of \eqref{subsec_build-adjoint}, the rule :
$j\mapsto\beta^\dagger_j$ yields a pseudo-natural transformation
$$
\beta^\dagger:G\Rightarrow 2\tdu\!\!\!\int_\phi F
$$
with coherence constraint given by the system of\/ $2$-cells
$\tau^{\beta^\dagger}_\bullet$.
\end{lemma}
\begin{proof} A simple inspection shows that
$\tau^{\beta^\dagger}_{\one_j}=\one_{\beta^\dagger_j}$ for every $j\in\Ob(J)$.
Next, let $g:j\to g'$ and $g':j'\to j''$ be two $1$-cells
of $J$; we have to check the identity :
$$
(\beta^\dagger_{j''}*\gamma^G_{g,g'})\odot(\tau^{\beta^\dagger}_{g'}*Gg)
\odot(L_{F,g'}*\tau^{\beta^\dagger}_g)=
\tau^{\beta^\dagger}_{g'\circ g}\odot(L_{F,g,g'}*\beta^\dagger_j).
$$
To this aim, set
$$
\begin{aligned}
X&\,:=(\Theta^{\beta,j}\circ g^*\circ g'^*)\odot
((\beta*\st_{j''})*(\Gamma^{G,g}\circ g'^*))\odot
((\beta*\st_{j''})*\Gamma^{G,g'}*\sF_{Gg}) \\
X'&\,:=(\Omega^{F,g'}*\sF_{L_{F,g}}*\sF_{\beta^\dagger_j})
\odot((\Omega^{F,g}\circ g'^*)*\sF_{\beta^\dagger_j})\odot X \\
Y&\,:=(\pi^{F,j''}*((\sF_{\beta^\dagger_{j''}}*\sF_{\gamma^G_{g,g'}})
\odot(\sF_{\tau^{\beta^\dagger}_{g'}}*\sF_{Gg})
\odot(\sF_{L_{F,g'}}*\sF_{\tau^{\beta^\dagger}_g})))\odot X' \\
Z&\,:=(\pi^{F,j''}*(\sF_{\tau^{\beta^\dagger}_{g'\circ g}}\odot
(\sF_{L_{F,g,g'}}*\sF_{\beta^\dagger_j})))
\odot X'.
\end{aligned}
$$
It suffices to show that $Y=Z$. Set as well
$Y':=\pi^{F,j''}*((\sF_{\beta^\dagger_{j''}}*\sF_{\gamma^G_{g,g'}})\odot
(\sF_{\tau^{\beta^\dagger}_{g'}}*\sF_{Gg}))$;
we compute :
$$
\begin{aligned}
Y&=Y'\!\odot\!(\Omega^{F,g'}*\sF_{\beta^\dagger_{j'}}*\sF_{Gg})
\!\odot\!((\pi^{F,j'}*g'^*)*\sF_{\tau^{\beta^\dagger}_g})
\odot((\Omega^{F,g}\circ g'^*)*\sF_{\beta^\dagger_j})\odot X' \\
&=Y'\odot(\Omega^{F,g'}*\sF_{\beta^\dagger_{j'}}*\sF_{Gg})
\odot((\Theta^{\beta,j'}\circ g'^*)*\sF_{Gg})\odot
((\beta*\st_{j''})*\Gamma^{G,g'}*\sF_{Gg}) \\
&=(\pi^{F,j''}*\sF_{\beta^\dagger_{j''}}*\sF_{\gamma^G_{g,g'}})\odot
(\Theta^{\beta,j''}*\sF_{Gg'}*\sF_{Gg}) \\
&=(\Theta^{\beta,j''}*\sF_{G(g'\circ g)})\odot
((\beta*\st_{j''})*(G*\phi^*_{j''})*\sF_{\gamma^G_{g,g'}})
\end{aligned}
$$
whereas : $Z=(\pi^{F,j''}*\sF_{\tau^{\beta^\dagger}_{g'\circ g}})\odot
(\Omega^{F,g'\circ g}*\sF_{\beta^\dagger_j})\odot X$. So, we are
reduced to checking that
$$
(\Gamma^{G,g}\circ g'^*)\odot
(\Gamma^{G,g'}*\sF_{Gg})=
\Gamma^{G,g'\circ g}\odot((G*\phi^*_{j''})*\sF_{\gamma^G_{g,g'}}).
$$
The latter follows from the composition axiom for $\gamma^G$.
Lastly, the naturality of $\tau^{\beta^\dagger}$ amounts to the
identity :
$$
\tau^{\beta^\dagger}_{g'}\odot(L_{F,\alpha}*\beta^\dagger_j)=
(\beta^\dagger_{j'}*G\alpha)\odot\tau^{\beta^\dagger}_g
\qquad
\text{for every $2$-cell $\alpha:g\Rightarrow g'$ of $J$.}
$$
To prove the latter, set $X:=(\Omega^{F,g}*\sF_{\beta^\dagger_j})
\odot(\Theta^{\beta,j}\circ g^*)\odot((\beta*\st_{j'})*\Gamma^{G,g})$;
it suffices to show:
$$
Y:=(\pi^{F,j'}*(\sF_{\tau^{\beta^\dagger}_{g'}}\odot
(\sF_{L_{F,\alpha}}*\sF_{\beta^\dagger_j})))\odot X=
Z:=(\pi^{F,j'}*((\sF_{\beta^\dagger_{j'}}*\sF_{G\alpha})\odot
\sF_{\tau^{\beta^\dagger}_g}))\odot X.
$$
To this aim, we remark :

\begin{claim}\label{cl_phi-unital-used-here}
$(\pi^{F,\alpha}*\sF_{\beta^\dagger_j})\odot(\Theta^{\beta,j}\circ g^*)=
(\Theta^{\beta,j}\circ g'^*)\odot((\beta*\st_{j'})*
(G\circ(\phi^*_j*\alpha^*))$.
\end{claim}
\begin{pfclaim} We have
$\alpha^*(f)=(\one_i,\alpha):g\circ f\to g'\circ f$ for every
$(f:j'\to\phi(i))\in\Ob(j'/\phi I)$. Now, let $\tau^{G*\phi^*_j}$
be the coherence constraint of $G*\phi^*_j$; since $\phi$ is
unital, it is easily seen that $\tau^{G*\phi^*_j}_{\alpha^*(f)}=G(f*\alpha)$
(details left to the reader). Then the compatibility condition
of $\Theta^{\beta,j}$ yields the commutative diagram :
$$
\xymatrix{ \beta_i\circ G(f\circ g)
\ar@{=>}[rr]^-{\Theta^{\beta,j}_{g^*(f)}} \ar@{=>}[d]_{\beta_i*G(f*\alpha)} & &
\pi^{F,j}_{g^*(f)}\circ\beta^\dagger_j \ar@{=>}[d]^{\pi^{F,\alpha}} \\
\beta_i\circ G(f\circ g') \ar@{=>}[rr]^-{\Theta^{\beta,j}_{g'^*(f)}} & &
\pi^{F,j}_{g'^*(f)}\odot\beta^\dagger_j 
}$$
whence the claim.
\end{pfclaim}

Using claim \ref{cl_phi-unital-used-here} we compute :
$$
\begin{aligned}
Y&\!=\!(\pi^{F,j'}*\sF_{\tau^{\beta^\dagger}_{g'}})\odot
((\Omega^{F,g'}\odot\pi^{F,\alpha})*\sF_{\beta^\dagger_j})\odot
(\Theta^{\beta,j}\circ g^*)\odot((\beta*\st_{j'})*\Gamma^{G,g}) \\
&\!=\!(\pi^{F,j'}*\sF_{\tau^{\beta^\dagger}_{g'}})\!\odot\!
(\Omega^{F,g'}*\sF_{\beta^\dagger_j})\!\odot\!(\Theta^{\beta,j}\circ g'^*)
\!\odot\!((\beta*\st_{j'})*(G\circ(\phi^*_j*\alpha^*))
\!\odot\!((\beta*\st_{j'})*\Gamma^{G,g}) \\
&\!=\!(\Theta^{\beta,j'}*\sF_{Gg'})\odot
((\beta*\st_{j'})*\Gamma^{G,g'})^{-1}\odot
((\beta*\st_{j'})*(G\circ(\phi^*_j*\alpha^*))\odot
((\beta*\st_{j'})*\Gamma^{G,g})
\end{aligned}
$$
and on the other hand
$$
Z=(\pi^{F,j'}*\sF_{\beta^\dagger_{j'}}*\sF_{G\alpha})\odot
(\Theta^{\beta,j'}\!\!*\!\sF_{Gg})=(\Theta^{\beta,j'}\!\!*\!\sF_{Gg'})
\odot((\beta*\st_{j'})*(G*\phi^*_{j'})*\sF_{G\alpha})
$$
so we are reduced to showing that :
$$
(G\circ(\phi^*_j*\alpha^*))\odot\Gamma^{G,g}=
\Gamma^{G,g'}\odot((G*\phi^*_{j'})*\sF_{G\alpha}).
$$
But the latter follows from the naturality condition for
$\gamma^G$.
\end{proof}

\begin{remark}\label{rem_cont-strong-2-Kan-ext}
In the situation of \eqref{subsec_build-adjoint}, suppose that
$\cC$ is strongly $2$-complete, so that we may take for
$2\tdu\!\!\int_\phi$ the strong $2$-right Kan extension of
remark \ref{rem_strong-2-Kan-ext}(ii). Then we may also choose
$\beta^\dagger_j$ such that $(\beta*\st_j)\odot(G*\phi^*_j)=
\pi^{F,j}\odot\sF_{\beta^\dagger_j}$, and we may let $\Theta^{\beta,j}$
be the corresponding identity modification. Suppose now
additionally that $G$ is a strict pseudo-functor; by inspecting
\eqref{eq_redo-sat-night}, we deduce that with such choices,
$\beta^\dagger$ is then a {\em strict} pseudo-natural transformation.
\end{remark}

\sset\subsubsection{}\label{subsec_build-adjoint-step-2}
Keep the notation of \eqref{subsec_build-adjoint}, and let
$\beta,\beta':G\circ\phi\Rightarrow F$ be two pseudo-natural
transformations and $\Xi:\beta\leadsto\beta'$ a modification.
Then there exists for every $j\in\Ob(J)$ a unique $2$-cell
$\Xi^\dagger_j:\beta^\dagger_j\Rightarrow\beta'^\dagger_j$ such that
\set\begin{equation}\label{eq_ouf-not-so-bad}
(\pi^{F,j}*\sF_{\Xi^\dagger_j})\odot\Theta^{\beta,j}=
\Theta^{\beta',j}\odot((\Xi\circ\st_j)*(G*\phi^*_j)).
\end{equation}

\begin{lemma}\label{lem_build-2-adj-step-2}
In the situation of \eqref{subsec_build-adjoint-step-2},
the rule : $j\mapsto\Xi^\dagger_j$ defines a modification
$$
\Xi^\dagger:\beta^\dagger\leadsto\beta'^\dagger.
$$
\end{lemma}
\begin{proof} Let $g:j\to j'$ be any $1$-cell of $J$; we
need to show the identity :
$$
\tau^{\beta'^\dagger}_g\odot(L_{F,g}*\Xi^\dagger_j)=
(\Xi^\dagger_{j'}*Gg)\odot\tau^{\beta^\dagger}_g.
$$
To this aim, set $X:=(\Omega^{F,g}*\sF_{\beta^\dagger_j})\odot
(\Theta^{\beta,j}\circ g^*)\odot((\beta*\st_{j'})*\Gamma^{G,g})$;
it suffices to show that
$$
Y:=(\pi^{F,j'}*(\sF_{\tau^{\beta'^\dagger}_g}\odot(\sF_{L_{F,g}}*\sF_{\Xi^\dagger_j)}))
\odot X=Z:=(\pi^{F,j'}*((\sF_{\Xi^\dagger_{j'}}*\sF_{Gg})\odot
\sF_{\tau^{\beta^\dagger}_g}))\odot X.
$$
We compute :
$$
\begin{aligned}
Y&=(\pi^{F,j'}*\sF_{\tau^{\beta'^\dagger}_g})\odot(\Omega^{F,g}*\sF_{\beta'^\dagger_j})
\odot((\pi^{F,j}*g^*)*\sF_{\Xi^\dagger_j})
\odot(\Theta^{\beta,j}\circ g^*)\odot((\beta*\st_{j'})*\Gamma^{G,g}) \\
&=(\pi^{F,j'}*\sF_{\tau^{\beta'^\dagger}_g})\!\odot\!(\Omega^{F,g}*\sF_{\beta'^\dagger_j})
\!\odot\!(\Theta^{\beta',j}\circ g^*)\!\odot\!
((\Xi\circ\st_{j'})*(G*\phi^*_j*g^*))\!\odot\!((\beta*\st_{j'})*\Gamma^{G,g}) \\
&=(\pi^{F,j'}\!*\!\sF_{\tau^{\beta'^\dagger}_g})\!\odot\!(\Omega^{F,g}*\sF_{\beta'^\dagger_j})
\!\odot\!(\Theta^{\beta',j}\!\circ\!g^*)\!\odot\!
((\beta'*\st_{j'})\!*\!\Gamma^{G,g})\!\odot\!
((\Xi\circ\st_{j'})*(G*\phi^*_{j'})\!*\!\sF_{Gg}) \\
&=(\Theta^{\beta',j'}*\sF_{Gg})\odot((\Xi\circ\st_{j'})*(G*\phi^*_{j'})*\sF_{Gg}).
\end{aligned}
$$
whereas : $Z=(\pi^{F,j'}*\sF_{\Xi^\dagger_{j'}}*\sF_{Gg})\odot
(\Theta^{\beta,j'}*\sF_{Gg})$, whence the contention.
\end{proof}

\begin{proposition}\label{prop_Kan-adjunction}
With the notation of lemmata {\em\ref{lem_build-2-adj-step-1}}
and {\em\ref{lem_build-2-adj-step-2}}, the rules :
$\beta\mapsto\beta^\dagger$ and $\Xi\mapsto\Xi^\dagger$
define a functor
$$
(-)^\dagger_{F,G}:\sPsNat(G\circ\phi,F)\to
\sPsNat\Bigl(G,2\tdu\!\!\!\int_\phi F\Bigr).
$$
\end{proposition}
\begin{proof} A simple inspection shows that
$(\one_\beta)^\dagger=\one_{\beta^\dagger}$ for every pseudo-natural
transformation $\beta:G\circ\phi\Rightarrow F$. It remains
to check that $(\Xi'\odot\Xi)^\dagger=\Xi'^\dagger\odot\Xi$
for every three pseudo-natural transformations
$\beta,\beta',\beta'':G\circ\phi\Rightarrow F$ and every pair
of modifications $\Xi:\beta\leadsto\beta'$ and
$\Xi':\beta'\leadsto\beta''$. However, we have :
$$
\begin{aligned}
(\pi^{F,j}*\sF_{\Xi'^\dagger_j}*\sF_{\Xi^\dagger_j})\odot\Theta^{\beta,j}
&=(\pi^{F,j}*\sF_{\Xi'^\dagger_j})\odot\Theta^{\beta',j}\odot
((\Xi\circ\st_j)*(G*\phi^*_j)) \\
&=\Theta^{\beta'',j}\odot((\Xi'\circ\st_j)*(G*\phi^*_j))\odot
((\Xi\circ\st_j)*(G*\phi^*_j)) \\
&=\Theta^{\beta'',j}\odot(((\Xi'\odot\Xi)\circ\st_j)*(G*\phi^*_j))
\end{aligned}
$$
whence the contention.
\end{proof}

\sset\subsubsection{}
Keep the notation of \eqref{subsec_build-adjoint}. We set
$$
\begin{aligned}
\omega^F_i&:=\pi^{F,\phi(i)}_{(i,\one_{\phi(i)})}:L_{F,\phi(i)}\to Fi
& &\text{for every $i\in\Ob(I)$} \\
\tau^{\omega^F}_f&:=\Omega^{F,\phi(f)}_{(i',\one_{\phi(i')})}\odot
\tau^{F,\phi(i)}_{(f,\one_{\phi(f)})}:
Ff\circ\omega^F_i\Rightarrow\omega^F_{i'}\circ L_{F,\phi(f)}
& &\text{for every $1$-cell $f:i\to i'$ of $I$}
\end{aligned}
$$
where $\tau^{F,\phi(i)}$ denotes the coherence constraint
of $\pi^{F,\phi(i)}$.

\begin{lemma}\label{lem_inverse-functor}
The rule : $i\mapsto\omega^F_i$ for every $i\in\Ob(I)$
defines a pseudo-natural transformation
$$
\omega^F:\Bigl(2\tdu\!\!\!\int_\phi F\Bigr)\circ\phi\Rightarrow F
$$
whose coherence constraint is given by the system of
$2$-cells $\tau^{\omega^F}_\bullet$.
\end{lemma}
\begin{proof} Let $f,f':i\to i'$ be two $1$-cells of $I$ and
$\alpha:f\Rightarrow f'$ a $2$-cell; for the naturality of
$\tau^{\omega^F}$ we need to check the identity :
$$
\tau^{\omega^F}_{f'}\odot(F\alpha*\omega^F_i)=
(\omega^F_{i'}*L_{F,\phi(\alpha)})\odot\tau^{\omega^F}_f.
$$
However, notice that $\alpha$ induces two $1$-cells
$(f,\phi(\alpha)),(f',\one_{\phi(f')}):(i,\one_{\phi(i)})\to(i',\phi(f'))$
and a $2$-cell $\alpha:(f,\phi(\alpha))\Rightarrow(f',\one_{\phi(f')})$
of $\phi(i)/\phi I$. By naturality of $\tau^{F,\phi(i)}$ we deduce the
identity :
$$
\tau^{F,\phi(i)}_{(f,\phi(\alpha))}=
\pi^{F,\phi(i)}_{(f',\one_{\phi(f')})}\odot(F\alpha*\pi^{F,\phi(i)}_{(i,\one_{\phi(i)})})
$$
whence :
$$
\tau^{\omega^F}_{f'}\odot(F\alpha*\omega^F_i)=
\Omega^{F,\phi(f')}_{(i',\one_{\phi(i')})}\odot\tau^{F,\phi(i)}_{(f,\phi(\alpha))}.
$$
On the other hand, from \eqref{eq_my-goodness} we get :
$$
(\omega^F_{i'}*L_{F,\phi(\alpha)})\odot\tau^{\omega^F}_f=
\Omega^{F,\phi(f')}_{(i',\one_{\phi(i')})}\odot
\pi^{F,\phi(\alpha)}_{(i',\one_{\phi(i')})}\odot
\tau^{F,\phi(i)}_{(f,\one_{\phi(f)})}.
$$
So we are reduced to checking the identity :
$\tau^{F,\phi(i)}_{(f,\phi(\alpha))}=\tau^{F,\phi(i)}_{(\one_{i'},\phi(\alpha))}
\odot\tau^{F,\phi(i)}_{(f,\one_{\phi(f)})}$. The latter follows from the
coherence axiom for $\tau^{F,\phi(i)}$. Next, we check the coherence
axioms for $\omega^F$. To this aim, notice that for every
$i\in\Ob(I)$ and $f:=\one_i$, the two-cell
$\tau^{F,\phi(i)}_{(f,\one_{\phi(f)})}$ is the identity of the $1$-cell
$\pi^{F,\phi(i)}_{(i,\one_{\phi(i)})}$ (remark \ref{rem_unital}(ii));
it follows easily that $\tau^{\omega^F}_{\one_i}=\one_{\omega^F_i}$.
Lastly, consider two $1$-cells $f:i\to i'$, $f':i'\to i''$ of
$I$, and set
$X:=\omega^F_{i''}*(L_{F,\gamma^\phi_{f,f'}}\odot L_{F,\phi(f),\phi(f')})$;
we need to show that
$$
Y:=X\odot(\tau^{\omega^F}_{f'}*L_{F,\phi(f)})\odot(Ff'*\tau^{\omega^F}_f)=
Z:=\tau^{\omega^F}_{f'\circ f}\odot(\gamma^F_{f,f'}*\omega^F_i)
$$
where $\gamma^\phi$ and $\gamma^F$ are the coherence constraints
of $\phi$ and $F$. We compute :
$$
\begin{aligned}
Y&=X\odot((\Omega^{F,\phi(f')}_{(i'',\one_{\phi(i'')})}\odot
\tau^{F,\phi(i')}_{(f',\one_{\phi(f')})})*L_{F,\phi(f)})\odot
(Ff'*(\Omega^{F,\phi(f)}_{(i',\one_{\phi(i')})}\odot
\tau^{F,\phi(i)}_{(f,\one_{\phi(f)})})) \\
&=X\odot(\Omega^{F,\phi(f')}_{(i'',\one_{\phi(i'')})}*L_{F,\phi(f)})\odot
\Omega^{F,\phi(f)}_{(i',\phi(f'))}\odot\tau^{F,\phi(i)}_{(f',\one_{\phi(f')\circ\phi(f)})}
\odot(Ff'*\tau^{F,\phi(i)}_{(f,\one_{\phi(f)})}) \\
&=(\omega^F_{i''}*L_{F,\gamma^\phi_{f,f'}})\odot
\Omega^{F,\phi(f')\circ\phi(f)}_{(i'',\one_{\phi(i'')})}
\odot\tau^{F,\phi(i)}_{(f',\one_{\phi(f')\circ\phi(f)})}
\odot(Ff'*\tau^{F,\phi(i)}_{(f,\one_{\phi(f)})}) \\
&=\Omega^{F,\phi(f'\circ f)}_{(i'',\one_{\phi(i'')})}\odot
\pi^{F,\gamma^\phi_{f,f'}}_{(i'',\one_{\phi(i'')})}\odot
\tau^{F,\phi(i)}_{(f',\one_{\phi(f')\circ\phi(f)})}
\odot(Ff'*\tau^{F,\phi(i)}_{(f,\one_{\phi(f)})}) \\
&=\Omega^{F,\phi(f'\circ f)}_{(i'',\one_{\phi(i'')})}\odot
\pi^{F,\gamma^\phi_{f,f'}}_{(i'',\one_{\phi(i'')})}\odot
\tau^{F,\phi(i)}_{(f'\circ f,(\gamma^\phi_{f,f'})^{-1})}\odot
(\gamma^F_{f,f'}*\pi^{F,\phi(i)}_{i,\one_{\phi(i)}})
\end{aligned}
$$
where the second equality follows from the compatibility
condition for $\Omega^{F,\phi(f)}$, the third follows from
\eqref{eq_triple-whammy}, the fourth follows from
\eqref{eq_my-goodness}, and the fifth holds by virtue of
the coherence axiom for $\tau^{F,\phi(i)}$, taking into
account \eqref{eq_def-composition}. So, we are reduced
to checking that
$$
\pi^{F,\gamma^\phi_{f,f'}}_{(i'',\one_{\phi(i'')})}\odot
\tau^{F,\phi(i)}_{(f'\circ f,(\gamma^\phi_{f,f'})^{-1}))}=
\tau^{F,\phi(i)}_{(f'\circ f,\one_{\phi(f'\circ f)})}
$$
which follows again from \eqref{eq_def-composition} and
the coherence axiom for $\tau^{F,\phi(i)}$.
\end{proof}

\sset\subsubsection{}\label{subsec_ddagger-functor}
In view of lemma \ref{lem_inverse-functor} we may define
a functor :
$$
(-)^\ddagger_{F,G}:\sPsNat\Bigl(G,2\tdu\!\!\!\int_\phi F\Bigr)\to
\sPsNat(G\circ\phi,F)
\qquad
\beta\mapsto\beta^\ddagger:=\omega^F\odot(\beta*\phi)
$$
attaching to every modification $\Xi:\beta\leadsto\beta'$
of pseudo-natural transformations
$\beta,\beta':G\Rightarrow 2\tdu\!\!\int_\phi F$, the
modification $\omega^F*(\Xi\circ\phi)$. For any
pseudo-natural transformation $\beta:G\circ\phi\Rightarrow F$
we set
$$
\Lambda^\beta_i:=\Theta^{\beta,\phi(i)}_{(i,\one_{\phi(i)})}:
\beta_i\Rightarrow(\beta^\dagger)^\ddagger_i
\qquad
\text{for every $i\in\Ob(I)$}.
$$

\begin{lemma}\label{lem_first-inverse}
The rule : $i\mapsto\Lambda^\beta_i$ defines an invertible
modification
$$
\Lambda^\beta:\beta\leadsto(\beta^\dagger)^\ddagger
$$
and the rule : $\beta\mapsto\Lambda^\beta$ defines an
isomorphism of functors :
$$
\Lambda:\one_{\sPsNat(G\circ\phi,F)}\isom
(-)^\ddagger_{F,G}\circ(-)^\dagger_{F,G}.
$$
\end{lemma}
\begin{proof} Let $f:i\to i'$ be any $1$-cell of $I$; we need
to check the identity :
$$
Y:=(\pi^{F,\phi(i')}_{(i',\one_{\phi(i')})}*\tau^{\beta^\dagger}_{\phi(f)})
\odot(\tau^{\omega^F}_f*\beta^\dagger_{\phi(i)})\odot
(Ff*\Theta^{\beta,\phi(i)}_{(i,\one_{\phi(i)})})=Z:=
(\Theta^{\beta,\phi(i')}_{i',\one_{\phi(i')}}*G\phi(f))\odot\tau^\beta_f
$$
where $\tau^\beta$ denotes the coherence constraint of $\beta$.
Let also $\tau^{(\beta*\st_{\phi(i)})\odot(G*\phi^*_{\phi(i)})}$ be the
coherence constraint of $(\beta*\st_{\phi(i)})\odot(G*\phi^*_{\phi(i)})$,
and notice that
$\tau^{(\beta*\st_{\phi(i)})\odot(G*\phi^*_{\phi(i)})}_{(f,\one_{\phi(f)})}=\tau^\beta_f$.
Combining with the compatibility condition for $\Theta^{\beta,\phi(i)}$
relative to the $1$-cell $(f,\one_{\phi(f)})$ of $\phi(i)/\phi I$,
we get :
$$
Y=(\pi^{F,\phi(i')}_{(i',\one_{\phi(i')})}*\tau^{\beta^\dagger}_{\phi(f)})
\odot(\Omega^{F,\phi(f)}_{(i',\one_{\phi(i')})}*\beta^\dagger_{\phi(i)})
\odot(\Theta^{\beta,\phi(i)}_{(i',\phi(f))}*G\phi(i))
\odot\tau^\beta_f.
$$
Then the sought identity follows from \eqref{eq_redo-sat-night},
applied with $g:=(f,\one_{\phi(i)})$, after noticing that
$\Gamma^{G,g}_{(i',\one_{\phi(i')})}=\one_{G\phi(f)}$, since $G$ is
unital. Lastly, let $\beta,\beta':G\Rightarrow 2\tdu\!\!\int_\phi F$
be two pseudo-natural transformations, and $\Xi:\beta\leadsto\beta'$
a modification; the naturality of $\Lambda$ amounts to the identity :
$$
(\omega^F_i*\Xi^\dagger_{\phi(i)})\odot\Lambda^\beta_i=
\Lambda^{\beta'}_i\odot\Xi_i
\qquad
\text{for every $i\in\Ob(I)$}
$$
which follows directly from \eqref{eq_ouf-not-so-bad}.
\end{proof}

\begin{lemma}
In the situation of \eqref{subsec_build-adjoint}, for every
$j\in\Ob(J)$ the rule
$$
(i,j\xrightarrow{g}\phi(i))\mapsto
\Delta^{F,j}_{(i,g)}:=\Omega^{F,g}_{(i,\one_{\phi(i)})}
\qquad
\text{for every $(i,g)\in\Ob(j/\phi I)$}
$$
defines an invertible modification
$$
\Delta^{F,j}:\pi^{F,j}\leadsto(\omega^F*\st_j)\odot
\Bigl(\Bigl(2\tdu\!\!\!\int_\phi F\Bigr)*\phi^*_j\Bigr).
$$
\end{lemma}
\begin{proof} Let $(f,\alpha):(i,g)\to(i',g')$ be any $1$-cell
of $j/\phi I$; we need to check the identity :
$$
\Omega^{F,g'}_{(i',\one_{\phi(i')})}\odot\tau^{F,j}_{(f,\alpha)}=
Y:=(\pi^{F,\phi(i')}_{(i',\one_{\phi(i')})}*(L_{F,\alpha}\odot L_{F,g,\phi(f)}))
\odot(\tau^{\omega^F}_f*L_{F,g})\odot(Ff*\Omega^{F,g}_{(i,\one_{\phi(i)})}).
$$
We compute :
$$
\begin{aligned}
Y&=(\pi^{F,\phi(i')}_{(i',\one_{\phi(i')})}*(L_{F,\alpha}\odot L_{F,g,\phi(f)}))
\odot((\Omega^{F,\phi(f)}_{(i',\one_{\phi(i')})}\odot
\tau^{F,\phi(i)}_{(f,\one_{\phi(f)})})*L_{F,g})\odot
(Ff*\Omega^{F,g}_{(i,\one_{\phi(i)})}) \\
&=(\pi^{F,\phi(i')}_{(i',\one_{\phi(i')})}*(L_{F,\alpha}\odot L_{F,g,\phi(f)}))
\odot(\Omega^{F,\phi(f)}_{(i',\one_{\phi(i')})}*L_{F,g})\odot
\Omega^{F,g}_{(i',\phi(f))}\odot\tau^{F,j}_{(f,g^*(\one_{\phi(f)}))} \\
&=(\pi^{F,\phi(i')}_{(i',\one_{\phi(i')})}*L_{F,\alpha})\odot
\Omega^{F,\phi(f)\circ g}_{(i',\one_{\phi(i')})}\odot
\tau^{F,j}_{(f,g^*(\one_{\phi(f)}))} \\
&=\Omega^{F,g'}_{(i',\one_{\phi(i')})}\odot\pi^{F,\alpha}_{(i',\one_{\phi(i')})}
\odot\tau^{F,j}_{(f,g^*(\one_{\phi(f)}))}
\end{aligned}
$$
where the second equality holds by the compatibility
condition of $\Omega^{F,g}$, applied to the $1$-cell
$(f,\one_{\phi(f)}):(i,\one_{\phi(i)})\to(i',\phi(f))$ of
$\phi(i)/\phi I$, the third equality follows from
\eqref{eq_triple-whammy}, and the fourth equality follows
from \eqref{eq_my-goodness}. So we are reduced to checking
that
$$
\tau^{F,j}_{(f,\alpha)}=\pi^{F,\alpha}_{(i',\one_{\phi(i')})}
\odot\tau^{F,j}_{(f,g^*(\one_{\phi(f)}))}
$$
which follows from the coherence condition for $\tau^{F,j}$,
relative to the composition of $1$-cells
$(i,g)\xrightarrow{(f,\one_{\phi(f)\circ g})}(i',\phi(f)\circ g)
\xrightarrow{(\one_{i'},\alpha)}(i',g')$.
\end{proof}

\sset\subsubsection{}
Now, let $\beta:G\Rightarrow 2\tdu\!\!\int_\phi F$ be a
pseudo-natural transformation. For every $j\in\Ob(J)$ we
get the $1$-cell $(\beta^\ddagger)^\dagger_j:Gj\to L_{F,j}$,
as well as the invertible modification
$$
\Theta^{\beta^\ddagger,j}:(\omega^F*\st_j)\odot(\beta*(\phi\circ\st_j))
\odot(G*\phi^*_j)\leadsto\pi^{F,j}\odot\sF_{(\beta^\ddagger)^\dagger_j}.
$$
On the other hand, by example \ref{ex_modifications}(ii)
we have the invertible modification
$$
\Upsilon^{\beta,j}:\Bigl(\Bigl(2\tdu\!\!\!\int_\phi F\Bigr)*\phi^*_j\Bigr)
\odot(\beta*\sF_j)\leadsto(\beta*(\phi\circ\st_j))\odot(G*\phi^*_j)
\qquad
(i,f:j\to\phi(i))\mapsto\tau^\beta_{\phi^*_{j,(i,f)}}
$$
where $\tau^\beta$ is the coherence constraint of $\beta$, and
notice that $\beta*\sF_j=\sF_{\beta_j}$. We set
$$
\Delta^{\beta,j}:=\Theta^{\beta^{\ddagger,j}}\odot
((\omega^F*\st_j)*\Upsilon^{\beta,j})\odot(\Delta^{F,j}*\sF_{\beta_j}):
\pi^{F,j}\odot\sF_{\beta_j}\leadsto\pi^{F,j}\odot\sF_{(\beta^\ddagger)^\dagger_j}.
$$
There exists then a unique invertible $2$-cell
$$
\Psi^\beta_j:\beta_j\Rightarrow(\beta^\ddagger)^\dagger_j
\qquad\text{such that}\qquad
\pi^{F,j}*\sF_{\Psi^\beta_j}=\Delta^{\beta,j}.
$$

\begin{lemma}\label{lem_second-inverse}
The rule : $j\mapsto\Psi^\beta_j$ for every $j\in\Ob(J)$
yields an invertible modification
$$
\Psi^\beta:\beta\leadsto(\beta^\ddagger)^\dagger.
$$
\end{lemma}
\begin{proof} Let $g:j\to j'$ be any $1$-cell of $J$;
we need to show that
$$
\tau^{(\beta^\ddagger)^\dagger}_g\odot(L_{F,g}*\Psi^\beta_j)=
(\Psi^\beta_{j'}*Gg)\odot\tau^\beta_g
$$
and recall that the coherence constraint $\tau^{(\beta^\ddagger)^\dagger}_g$
of $(\beta^\ddagger)^\dagger$ is characterized by the identity :
$$
(\pi^{F,j'}*\sF_{\tau^{(\beta^\ddagger)^\dagger}_g})\odot
(\Omega^{F,g}*\sF_{(\beta^\ddagger)^\dagger_j})\odot
(\Theta^{\beta^\ddagger,j}\circ g^*)\odot
((\beta^\ddagger*\st_{j'})*\Gamma^{G,g})
=\Theta^{\beta^\ddagger,j'}*\sF_{Gg}.
$$
As usual, we reduce to checking the identity :
$$
Y:=\pi^{F,j'}*(\sF_{\tau^{(\beta^\ddagger)^\dagger}_g}\odot
(\sF_{L_{F,g}}*\sF_{\Psi^\beta_j}))=
Z:=\pi^{F,j'}*((\sF_{\Psi^\beta_{j'}}*\sF_{Gg})\odot\sF_{\tau^\beta_g}).
$$
To ease notation, set $H:=2\tdu\!\!\int_\phi F$; we remark :

\begin{claim}\label{cl_from-coh-tau-beta}
$((\beta*\phi*\st_{j'})*\Gamma^{G,g})^{-1}\odot
(\Upsilon^{\beta,j}\circ g^*)=(\Upsilon^{\beta,j'}*\sF_{Gg})\odot
((H*\phi^*_{j'})*\sF_{\tau^\beta_g})\odot(\Gamma^{H,g}*\sF_{\beta_j})^{-1}$.
\end{claim}
\begin{pfclaim} The left-hand side is the modification
$$
(H*\phi^*_j*g^*)\odot\sF_{\beta_j}\leadsto
(\beta*(\phi\circ\st_{j'}))\odot(G*\phi^*_{j'})\odot\sF_{Gg}
\qquad
(i,h:j'\to\phi(i))\mapsto\tau^\beta_{h\circ g}\odot
(\beta_{\phi(i)}*\gamma^G_{g,h})^{-1}
$$
and the coherence axiom for $\tau^\beta$ gives :
$$
(\beta_{\phi(i)}*\gamma^G_{g,h})^{-1}\odot\tau^\beta_{h\circ g}=
(\tau^\beta_h*Gg)\odot(L_{F,h}*\tau^\beta_g)\odot(L_{F,g,h}*\beta_j)^{-1}
$$
where $\gamma^G$ is the coherence constraint of $G$. The claim
translates the latter identity.
\end{pfclaim}

Now, let $X:=\Theta^{\beta^\ddagger,j'}*\sF_{Gg}$ and
$X':=X\odot(\omega^F*\st_{j'})*\Upsilon^{\beta',j}*\sF_{Gg}$. We compute :
$$
\begin{aligned}
Y&=X\odot((\beta^\ddagger*\st_{j'})*\Gamma^{G,g})^{-1}\odot
(\Theta^{\beta^\ddagger,j}\circ g^*)^{-1}\odot
(\Omega^{F,g}*\sF_{(\beta^\ddagger)^\dagger_j})^{-1}\odot
(\pi^{F,j'}*\sF_{L_{F,g}}*\sF_{\Psi^\beta_j}) \\
&=X\odot((\beta^\ddagger*\st_{j'})*\Gamma^{G,g})^{-1}\odot
(\Theta^{\beta^\ddagger,j}\circ g^*)^{-1}\odot
((\pi^{F,j}*g^*)*\sF_{\psi^\beta_j})\odot(\Omega^{F,g}*\sF_{\beta_j})^{-1} \\
&=X\odot((\beta^\ddagger*\st_{j'})*\Gamma^{G,g})^{-1}\odot
(\Theta^{\beta^\ddagger,j}\circ g^*)^{-1}\odot
(\Delta^{\beta,j}\circ g^*)\odot(\Omega^{F,g}*\sF_{\beta_j})^{-1} \\
&=X\!\odot\!((\beta^\ddagger*\st_{j'})*\Gamma^{G,g})^{-1}\!\odot\!
((\omega^F*\st_{j'})*(\Upsilon^{\beta,j}\circ g^*))\!\odot\!
((\Delta^{F,j}\circ g^*)*\sF_{\beta_j})\!\odot\!
(\Omega^{F,g}*\sF_{\beta_j})^{-1} \\
&=X'\!\odot\!((\omega^F\!\!*\st_{j'})\!*\!
(((H\!*\!\phi^*_{j'})*\sF_{\tau^\beta_g})\!\odot\!
(\Gamma^{H,g}*\sF_{\beta_j})^{-1}))\!\odot\!
((\Delta^{F,j}\!\circ\!g^*)*\sF_{\beta_j})\!\odot\!
(\Omega^{F,g}*\sF_{\beta_j})^{-1}
\end{aligned}
$$
where the last equality follows from claim
\ref{cl_from-coh-tau-beta}. On the other hand :
$$
Z=X'\odot(\Delta^{F,j'}*\sF_{\beta_{j'}}*\sF_{Gg})\odot
(\pi^{F,j'}*\sF_{\tau^\beta_g})=X'\odot((\omega^F*\st_{j'})*
(H*\phi^*_{j'})*\sF_{\tau^\beta_g})\odot(\Delta^{F,j'}*\sF_{Hg}*\sF_{\beta_j})
$$
so we are reduced to showing that
$$
\Delta^{F,j}\circ g^*=((\omega^F*\st_{j'})*\Gamma^{H,g})\odot
(\Delta^{F,j'}*\sF_{Hg})\odot\Omega^{F,g}
$$
which follows from \eqref{eq_triple-whammy}.
\end{proof}

\begin{proposition} For every unital pseudo-functors
$F:I\to\cC$ and $G:J\to\cC$, the functors $(-)^\dagger_{F,G}$
of proposition {\em\ref{prop_Kan-adjunction}} and $(-)^\ddagger_{F,G}$
of \eqref{subsec_ddagger-functor} are equivalences.
\end{proposition}
\begin{proof} In view of lemmata \ref{lem_first-inverse} and
\ref{lem_second-inverse}, it suffices to prove that the
rule : $\beta\mapsto\Psi^\beta$ defines a natural transformation
$$
\Psi:\one_{\sPsNat(G,2\tdu\!\!\int_\phi F)}\Rightarrow
(-)^\dagger_{F,G}\circ(-)^\ddagger_{F,G}.
$$
Thus, let $\alpha,\beta:G\Rightarrow 2\tdu\!\!\int_\phi F$ be two
pseudo-natural transformations, and $\Xi:\alpha\leadsto\beta$
a modification; we need to check the identity :
$$
\Psi^\beta_j\odot\Xi_j=
(\omega^F*(\Xi\circ\phi))^\dagger_j\odot\Psi^\alpha_j
\qquad
\text{for every $j\in\Ob(J)$}
$$
and as usual, it suffices to show that
$$
Y:=\pi^{F,j}*(\sF_{\Psi^\beta_j}\odot\sF_{\Xi_j})=
Z:=\pi^{F,j}*(\sF_{(\omega^F*(\Xi\circ\phi))^\dagger_j}\odot\sF_{\Psi^\alpha_j}).
$$
To ease notation, set $H:=2\tdu\!\!\int_\phi F$; we remark :

\begin{claim}\label{cl_poor-Pat}
$((\Xi\circ\phi\circ\st_j)*(G*\phi^*_j))\odot\Upsilon^{\alpha,j}=
\Upsilon^{\beta,j}\odot((H*\phi^*_j)*(\Xi\circ\sF_j))$.
\end{claim}
\begin{pfclaim} The left-hand side is the modification :
$$
(H*\phi^*_j)\odot(\alpha*\sF_j)\leadsto
(\beta*(\phi\circ\st_j))\odot(G*\phi^*_j)
\qquad
(j,g\mapsto\phi(i))\mapsto(\Xi_{\phi(i)}*Gg)\odot\tau^\alpha_g
$$
where $\tau^\alpha$ is the coherence constraint of $\alpha$;
then the compatibility condition for $\Xi$ gives :
$$
(\Xi_{\phi(i)}*Gg)\odot\tau^\alpha_g=\tau^\beta_g\odot(Hg*\Xi_j)
$$
and the claim translates the latter identity.
\end{pfclaim}

We have : $Y\!=\!\Delta^{\beta,j}\!\odot\!(\pi^{F,j}*\sF_{\Xi_j})\!=\!
\Theta^{\beta^{\ddagger,j}}\!\odot\!((\omega^F*\st_j)*\Upsilon^{\beta,j})
\!\odot\!(\Delta^{F,j}*\sF_{\beta_j})\!\odot\!(\pi^{F,j}*\sF_{\Xi_j})$, and
$$
\begin{aligned}
Z&=(\pi^{F,j}*\sF_{(\omega^F*(\Xi\circ\phi))^\dagger_j})\odot\Delta^\alpha_j \\
&=(\pi^{F,j}*\sF_{(\omega^F*(\Xi\circ\phi))^\dagger_j})\odot
\Theta^{\alpha^{\ddagger,j}}\odot
((\omega^F*\st_j)*\Upsilon^{\alpha,j})\odot(\Delta^{F,j}*\sF_{\alpha_j}) \\
&=\Theta^{\beta^\ddagger,j}\odot
(((\omega^F*(\Xi\circ\phi))\circ\st_j)*(G*\phi^*_j))\odot
((\omega^F*\st_j)*\Upsilon^{\alpha,j})\odot(\Delta^{F,j}*\sF_{\alpha_j}) \\
&=\Theta^{\beta^\ddagger,j}\odot
((\omega^F*\st_j)*(\Xi\circ\phi\circ\st_j)*(G*\phi^*_j))\odot
((\omega^F*\st_j)*\Upsilon^{\alpha,j})\odot(\Delta^{F,j}*\sF_{\alpha_j}) \\
&=\Theta^{\beta^\ddagger,j}\odot
((\omega^F*\st_j)*(\Upsilon^{\beta,j}\odot((H*\phi^*_j)*(\Xi\circ\sF_j))))
\odot(\Delta^{F,j}*\sF_{\alpha_j})
\end{aligned}
$$
where the last equality holds by claim \ref{cl_poor-Pat}. Thus,
we are reduced to checking the identity :
$$
\Delta^{F,j}*\sF_{\beta_j}\odot(\pi^{F,j}*\sF_{\Xi_j})=
((\omega^F*\st_j)*(H*\phi^*_j)*(\Xi\circ\sF_j))
\odot(\Delta^{F,j}*\sF_{\alpha_j})
$$
which follows as usual from remark \ref{rem_equiv-2-cat}(i).
\end{proof}

\begin{theorem}\label{th_2-Kan-extensions}
Let $I$ be a small $2$-category, $J$ a $2$-category with small
$\Hom$-categories, $\phi:I\to J$ any pseudo-functor, and suppose
that all the $2$-cells of $J$ are invertible. Then for every
$2$-complete (resp. $2$-cocomplete) $2$-category $\cC$ the right
(resp. left) $2$-Kan extension along $\phi$ is a right (resp.
left) $2$-adjoint for the pseudo-functor
$$
\sPsFun(\phi,\cC):\sPsFun(J,\cC)\to\sPsFun(I,\cC).
$$
\end{theorem}
\begin{proof} Let us show first that the rule :
$(F,G)\mapsto(-)^\ddagger_{F,G}$ yields a pseudo-natural equivalence :
$$
(-)^\ddagger:\sPsNat\Bigl(\one_{\uniPsFun(J,\cC)},2\tdu\!\!\!\int_\phi\Bigr)
\isom\sPsNat(\uniPsFun(\phi^u,\cC),\one_{\uniPsFun(I,\cC)})
$$
(notation of remark \ref{rem_uniPsFun}(i)). To this aim let
us remark :

\begin{claim} Let $F,F':I\to\cC$ be two unital pseudo-functors.
Every pseudo-natural transformation $\mu:F\Rightarrow F'$ induces
an invertible modification :
$$
\Sigma^\mu:\mu\odot\omega^F\leadsto
\omega^{F'}\odot\Bigl(\Bigl(2\tdu\!\!\!\int_\phi\mu\Bigr)*\phi\Bigr)
\qquad
i\mapsto\Omega^{\mu,\phi(i)}_{(i,\one_{\phi(i)})}.
$$
\end{claim}
\begin{pfclaim} Let $f:i\to i'$ be any $1$-cell of $I$,
and $\tau^\mu$ the coherence constraint of $\mu$; set :
$$
\begin{aligned}
Y&:=(\pi^{F',\phi(i')}_{(i',\one_{\phi(i')})}*L_{\mu,\phi(f)})\odot
(\tau^{\omega^{F'}}_f*L_{\mu,\phi(i)})\odot
(F'f*\Omega^{\mu,\phi(i)}_{(i,\one_{\phi(i)})}) \\
Z&:=(\Omega^{\mu,\phi(i')}_{(i',\one_{\phi(i')})}*L_{F,\phi(f)})
\odot(\mu_{i'}*\tau^{\omega^F}_f)\odot
(\tau^\mu_f*\pi^{F,\phi(i)}_{(i,\one_{\phi(i)})})
\end{aligned}
$$
so that we need to check the identity $Y=Z$. We compute :
$$
\begin{aligned}
Y&=(\pi^{F',\phi(i')}_{(i',\one_{\phi(i')})}*L_{\mu,\phi(f)})\odot
((\Omega^{F',\phi(f)}_{(i',\one_{\phi(i')})}\odot
\tau^{F',\phi(i)}_{(f,\one_{\phi(f)})})*L_{\mu,\phi(i)})\odot
(F'f*\Omega^{\mu,\phi(i)}_{(i,\one_{\phi(i)})}) \\
&=(\pi^{F',\phi(i')}_{(i',\one_{\phi(i')})}*L_{\mu,\phi(f)})\!\odot\!
(\Omega^{F',\phi(f)}_{(i',\one_{\phi(i')})}*L_{\mu,\phi(i)})\!\odot\!
\Omega^{\mu,\phi(i)}_{(i',\phi(f))}\!\odot\!
(\mu_{i'}*\tau^{F,\phi(i)}_{(f,\one_{\phi(f)})})\!\odot\!
(\tau^\mu_f*\pi^{F,\phi(i)}_{(i,\one_{\phi(i)})})
\end{aligned}
$$
where the second equality follows from the compatibility
condition of $\Omega^{\mu,\phi(i)}$, relative to the $1$-cell
$(f,\one_{\phi(f)}):(i,\one_{\phi(i)})\to(i',\phi(f))$ of
$\phi(i)/\phi I$. We are thus reduced to showing that :
$$
(\pi^{F',\phi(i')}_{(i',\one_{\phi(i')})}*L_{\mu,\phi(f)})\!\odot\!
(\Omega^{F',\phi(f)}_{(i',\one_{\phi(i')})}*L_{\mu,\phi(i)})\!\odot\!
\Omega^{\mu,\phi(i)}_{(i',\phi(f))}=
(\Omega^{\mu,\phi(i')}_{(i',\one_{\phi(i')})}*L_{F,\phi(f)})\!\odot\!
(\mu_{i'}*\Omega^{F,\phi(f)}_{(i',\one_{\phi(i')})})
$$
which follows from \eqref{eq_laughing-hard}.
\end{pfclaim}

\begin{claim} Let $F,F':I\to\cC$ and $G,G':J\to\cC$ be four unital
pseudo-functors, $\lambda:G'\Rightarrow G$ and $\mu:F\Rightarrow F'$
two pseudo-natural transformations. The rule
$$
\beta\mapsto\tau^\ddagger_{(\lambda,\mu),\beta}:=
\Sigma^\mu*((\beta\odot\lambda)*\phi)
$$
yields the orientation for an essentially commutative square
diagram of functors :
$$
\xymatrix{ \sPsNat(G,2\tdu\!\!\int_\phi F)
\ar[d]_{\sPsNat(\lambda,2\tdu\!\!\int_\phi\mu)}
\ar[rrr]^-{(-)^\ddagger_{F,G}} &
\drtwocell\omit{\qquad\tau^\ddagger_{(\lambda,\mu)}}
& &\sPsNat(G\circ\phi,F) \ar[d]^{\sPsNat(\lambda*\phi,\mu)} \\
\sPsNat(G',2\tdu\!\!\int_\phi F') \ar[rrr]_-{(-)^\ddagger_{F',G'}}
& & & \sPsNat(G'\circ\phi,F').
}$$
\end{claim}
\begin{pfclaim} Let $\beta,\beta':G\to 2\tdu\!\!\int_\phi F$ be
two pseudo-natural transformations, and $\Xi:\beta\leadsto\beta'$
a modification; we need to check the identity :
$$
\tau^\ddagger_{(\lambda,\mu),\beta'}\odot
(\mu*\omega^F*(\Xi\circ\phi)*(\lambda*\phi))=
\Bigl(\omega^{F'}*\Bigl(\Bigl(\Bigl(2\tdu\!\!\!\int_\phi\mu\Bigr)*
\Xi*\lambda\Bigr)\circ\phi\Bigr)\Bigr)
\odot\tau^\ddagger_{(\lambda,\mu),\beta}
$$
and we come down to showing that :
$$
(\Sigma^\mu*(\beta'*\phi))\odot(\mu*\omega^F*(\Xi\circ\phi))=
\Bigl(\omega^{F'}*\Bigl(\Bigr(2\tdu\!\!\!\int_\phi\mu\Bigr)*\phi\Bigr)
*(\Xi\circ\phi)\Bigr)\odot(\Sigma^\mu*(\beta*\phi))
$$
which follows from remark \ref{rem_equiv-2-cat}(i).
\end{pfclaim}

To conclude the proof, it remains to check that the rule
$(\lambda,\mu)\mapsto\tau^\ddagger_{(\lambda,\mu)}$ yields a
coherence constraint for the sought pseudo-functor $(-)^\ddagger$.
For the first coherence axiom, notice that both the source
and target of $(-)^\ddagger$ are themselves unital pseudo-functors,
hence it suffices to check that $\tau^\ddagger_{(\one_G,\one_F)}$
is the identity automorphism of $(-)_{F,G}^\ddagger$, for every
pair of unital pseudo-functors $F:I\to\cC$ and $G:J\to\cC$
(remark \ref{rem_unital}(ii)). We then come down to checking
that $\Sigma^{\one_F}=\one_{\omega^F}$ for every such $F$; the
latter identity follows by a direct inspection of the definitions.

Next, consider unital pseudo-functors $F,F',F'':I\to\cC$,
$G,G',G'':J\to\cC$ and two pairs of pseudo-natural transformations
$(\lambda':G''\Rightarrow G',\mu':F'\Rightarrow F'')$ and
$(\lambda:G'\Rightarrow G,\mu:F\Rightarrow F')$; we need to
check the identity :
$$
((\!-\!)^\ddagger_{F'',G''}*\sPsNat(G''\!,\!L_{\mu,\mu'}))\!\odot\!
\Bigl(\tau^\ddagger_{(\lambda',\mu')}*
\sPsNat\Bigl(\!\lambda,2\tdu\!\!\!\int_\phi\!\!\mu\!\Bigr)\!\Bigr)
\!\odot\!(\sPsNat(\lambda'*\phi,\mu')*\tau^\ddagger_{(\lambda,\mu)})
\!=\!\tau^\ddagger_{(\lambda'\odot\lambda,\mu'\odot\mu)}.
$$
Thus, let $\beta:G\Rightarrow 2\tdu\!\!\int_\phi F$ be
any pseudo-natural transformation, and set
$\rho:=(\beta\odot\lambda\odot\lambda')*\phi$; the assertion
comes down to the identity :
$$
(\omega^F*(L_{\mu,\mu'}\circ\phi)*\rho)
\odot\Bigl(\Sigma^{\mu'}*
\Bigl(\Bigl(\Bigl(2\tdu\!\!\!\int_\phi\mu\Bigr)*\phi\Bigr)
\odot\rho\Bigr)\Bigl)\odot(\mu'*\Sigma^\mu*\rho)=
\Sigma^{\mu'\odot\mu}*\rho
$$
so it suffices to show that :
$$
(\omega^F*(L_{\mu,\mu'}\circ\phi))
\odot\Bigl(\Sigma^{\mu'}*
\Bigl(\Bigl(2\tdu\!\!\!\int_\phi\mu\Bigr)*\phi\Bigr)\Bigl)
\odot(\mu'*\Sigma^\mu)=\Sigma^{\mu'\odot\mu}
$$
which holds by \eqref{eq_L_beta-beta-prime-j}. Lastly,
the source and target of $(-)^\ddagger$ are pseudo-functors
$$
\uniPsFun(J,\cC)^o\times\uniPsFun(I,\cC)\to\bCat.
$$
On the other hand, the strict $2$-equivalences $(-)^u$ of
remark \ref{rem_uniPsFun}(i) yield a strict $2$-equivalence
$$
U:\sPsFun(J,\cC)^o\times\sPsFun(I,\cC)\isom
\uniPsFun(J,\cC)^o\times\uniPsFun(I,\cC)
$$
and the strict isomorphisms of pseudo-functors
\eqref{eq_really-strict} yield pseudo-natural isomorphisms
$$
\begin{aligned}
\sPsNat\Bigl(\one_{\uniPsFun(J,\cC)},2\tdu\!\!\!\int_\phi\Bigr)\circ U
&\,\isom
\sPsNat\Bigl(\one_{\sPsFun(J,\cC)},2\tdu\!\!\!\int_\phi\Bigr) \\
\sPsNat(\uniPsFun(\phi^u,\cC),\one_{\uniPsFun(I,\cC)})\circ U
&\,\isom
\sPsNat(\sPsFun(\phi,\cC),\one_{\sPsFun(I,\cC)}).
\end{aligned}
$$
Summing up, we deduce the sought $2$-adjunction between
$\sPsFun(\phi,\cC)$ and $2\tdu\!\!\int_\phi$.

The assertion concerning the left $2$-Kan extension along
$\phi$ is an immediate consequence, in view of remark
\ref{rem_opposing-thumb}(i).
\end{proof}

\begin{corollary}\label{cor_2-Kan-ext-of-fully-faith}
In the situation of theorem {\em\ref{th_2-Kan-extensions}},
suppose that $\phi$ is fully faithful. Then the same holds
for the pseudo-functor $2\tdu\!\!\int_\phi$ (resp. for
$2\tdu\!\!\int^\phi$).
\end{corollary}
\begin{proof} As usual, it suffices to check the assertion
concerning right $2$-Kan extensions. To this aim, we notice :

\begin{claim}\label{cl_found-pseudo-initial}
Suppose that $\phi:I\to J$ is fully faithful and that all the
$2$-cells of $J$ are invertible. Then for every $i\in\Ob(I)$,
the object $(i,\one_{\phi(i)})$ is pseudo-initial in $\phi(i)/\phi I$.
\end{claim}
\begin{pfclaim} Let $(i',f)$ be any othe object of
$\phi(i)/\phi I$; by assumption, there exist a $1$-cell
$g:i\to i'$ in $I$ and an invertible $2$-cell
$\alpha:\phi(g)\isom f$ in $J$. Then we get the $1$-cell
$$
(i,\one_{\phi(i)})\xrightarrow{(g,\one_{\phi(g)})}(i',\phi(g))
\xrightarrow{(\one_{i'},\alpha)}(i',f).
$$
Next, let $(t,\alpha),(s,\beta):(i,\one_{\phi(i)})\to(i',f)$
be two $1$-cells in $\phi(i)/\phi I$; by assumption $\beta$
and $\alpha$ are invertible, hence we get the invertible
$2$-cell $\beta^{-1}\odot\alpha:\phi(t)\isom\phi(s)$, and
since $\phi$ is fully faithful, there exists a unique
invertible $2$-cell $\lambda:t\isom s$ in $I$ such that
$\phi(\lambda)=\beta^{-1}\odot\alpha$. We deduce an invertible
$2$-cell $\lambda:(t,\alpha)\isom(s,\beta)$ in $\phi(i)/\phi I$,
and it is easily seen that this is the unique $2$-cell from
$(t,\alpha)$ to $(s,\beta)$. The claim follows.
\end{pfclaim}

Let now $F:I\to\cC$ be any pseudo-functor, and for every
$j\in\Ob(J)$, let $(L_{F,j},\pi^{F,j})$ be a $2$-limit of
$F\circ\st_j$, as in \eqref{sec_2-cat-Kan-ext}; in view
of claim \ref{cl_found-pseudo-initial} and proposition
\ref{prop_pseudo-initial}, we deduce that the $1$-cell
$\pi^{F,\phi(i)}_{(i,\one_{\phi(i)})}:L_{F,j}\to Fi$ is an equivalence
for every $i\in\Ob(I)$, and therefore the pseudo-natural
transformation $\omega^F$ of lemma \ref{lem_inverse-functor}
is a pseudo-natural equivalence. To conclude, it suffices now
to invoke corollary \ref{cor_fully-faith-2-adjoint}.
\end{proof}

\section{Special categories}

\subsection{Fibrations}
\label{sec_phi-cartesian}
Let $\phi:\cA\to\cB$ be a functor, $f:A'\to A$ a morphism
in $\cA$. We say that $f$ is {\em $\phi$-cartesian}, or --
slightly abusively -- that $f$ is {\em $\cB$-cartesian\/},
if the induced commutative diagram of sets (notation of
\eqref{eq_push-for}) :
$$
\xymatrix{
\Hom_\cA(X,A') \ar[rr]^-{f_*} \ar[d]_\phi & &
\Hom_\cA(X,A) \ar[d]^\phi \\
\Hom_\cB(\phi X,\phi A') \ar[rr]^-{(\phi f)_*} & &
\Hom_\cB(\phi X,\phi A)
}$$
is cartesian for every $X\in\Ob(\cA)$. In this case, one
also says that {\em $f$ is an inverse image of $A$\/} over
$\phi f$, or -- slightly abusively -- that {\em $A'$ is
an inverse image of $A$\/} over $\phi f$. 

\begin{remark}\label{rem_first-rem-fibrations}
Keep the notation of \eqref{sec_phi-cartesian}.

(i)\ \ 
It is easily seen that the composition of two
$\cB$-cartesian morphisms is $\cB$-cartesian.

(ii)\ \
Let $g:B'\to B$ and $g':B''\to B'$ be two morphisms of
$\cB$, and $f:A'\to A$ and $f'':A''\to A$ two $\cB$-cartesian
morphisms of $\cA$ such that $\phi f=g$ and $\phi f''=g\circ g'$.
Then there exists a unique morphism $f':A''\to A'$ of $\cA$
such that $\phi f'=g'$ and $f''=f\circ f'$, and this morphism
is $\cB$-cartesian : the detailed verification shall be left
to the reader. 
\end{remark}

\begin{definition}\label{def_pull-bekko}
Let $\phi:\cA\to\cB$ be a functor, and $B$ any object of $\cB$.
\begin{enumerate}
\item
The {\em fibre of\/ $\phi$ over $B$\/} is the category whose
objects are all the $A\in\Ob(\cA)$ such that $\phi A=B$, and
whose morphisms $f:A'\to A$ are the elements of $\Hom_\cA(A',A)$
such that $\phi f=\one_B$. We denote this category by
$$
\phi^{-1}B
\qquad\text{or more simply by}\qquad
\cA_B
$$
in the contexts where the latter notation does not give
rise to ambiguities. We denote the natural faithful embedding
of $\cA_B$ into $\cA$ by :
$$
\iota_B:\cA_B\to\cA.
$$
\item
We say that $\phi$ is a {\em fibration\/} if, for every
morphism $g:B'\to B$ in $\cB$, and every $A\in\Ob(\cA_B)$,
there exists an inverse image $f:A'\to A$ of $A$ over $g$.
In this case, we also say that $\cA$ is a {\em fibred
$\cB$-category}, and $\phi$ is called the {\em structure
functor} of $\cA$.
\end{enumerate}
\end{definition}

\begin{example}\label{ex_fibred-cats}
Let $\cB$, $\cC$ be any two categories, $F:\cC\to\cB$ a functor,
$X$ any object of $\cB$.

(i)\ \
The source functor $\ss_X:F\cC/X\to\cC$ of \eqref{subsec_fibreovercat}
is a fibration, and all the morphisms in $F\cC/X$ are $\cC$-cartesian.
The easy verification shall be left to the reader.

(ii)\ \
The source functor $\ss:\cB/F\cC\to\cB$ of \eqref{subsec_target-fctr}
is a fibration. Namely, the $\ss$-cartesian morphisms are the
commutative diagrams
$$
\xymatrix{ B \ar[r]^-f \ar[d]_g & FC \ar[d]^{Fg'} \\
B' \ar[r]^-{f'} & FC'
}$$
such that $g'$ is an isomorphism in $\cC$. The fibre $\ss^{-1}X$
is the category $X/F\cC$. As a special case, the source functor
$\ss:\sMorph(\cB)\to\cB$ of \eqref{subsec_Morph-cat} is a fibration.

(iii)\ \
Suppose that all fibre products are representable in $\cB$. Then
also the target functor $\st:\sMorph(\cB)\to\cB$ is a fibration;
more precisely, the $\st$-cartesian morphisms are the square
diagrams as in (ii) (with $F:=\one_\cB$) which are cartesian
({\em i.e.} fibred). We have $\st^{-1}X=\cB/X$.
\end{example}

\begin{definition}\label{def_decartes-functors}
Let $\cA$, $\cA'$, and $\cB$ be three categories, $\phi:\cA\to\cB$
and $\phi':\cA'\to\cB$ two functors, and $F:\cA\to\cA'$ a
{\em $\cB$-functor}, {\em i.e.} $F$ verifies the identity
$\phi'\circ F=\phi$.

(i)\ \
We say that $F$ is {\em cartesian} if it sends $\cB$-cartesian
morphisms of $\cA$ to $\cB$-cartesian morphisms in $\cA'$. We
denote by :
$$
\sCart_\cB(\cA,\cA')
$$
the category whose objects are the cartesian $\cB$-functors
$F:\cA\to\cA'$, and whose morphisms are the {\em natural
$\cB$-transformations}, {\em i.e.} the natural transformations
such that :
$$
\alpha:F\Rightarrow G
\qquad\text{such that}\qquad
\phi'*\alpha=\one_\phi.
$$
The composition law is the usual composition of natural
transformations : $(\beta,\alpha)\mapsto\alpha\odot\beta$.
Notice that if $\cA$ and $\cA'$ are small, the same holds
for $\sCart_\cB(\cA,\cA')$ : details left to the reader.

(ii)\ \
For any two other $\cB$-categories $\cC\to\cB$ and
$\cC'\to\cB$, every pair of cartesian $\cB$-functors
$H:\cC\to\cA$ and $K:\cA'\to\cC'$ induces a functor :
$$
\sCart_\cB(H,K):\sCart_\cA(\cA,\cA')\to\sCart_\cB(\cC,\cC')
\qquad
G\mapsto K\circ G\circ H.
$$
To any morphism $\alpha:G\Rightarrow G'$ in
$\sCart_\cB(\cA,\cA')$, the functor $\sCart_\cB(H,K)$
assigns the natural transformation
$K*\alpha*H:K\circ G\circ H\Rightarrow K\circ G'\circ H$.
In case $H=\one_\cA$ (resp. $K=\one_{\cA'}$) we also denote
this functor by $\sCart_\cB(\cA,K)$ (resp. by $\sCart_\cB(H,\cA')$).

Likewise, if $H,H':\cC\to\cA$ and $K,K':\cA'\to\cC'$ are
four $\cB$-cartesian functors, every pair of natural
$\cB$-transformations $\beta:H\Rightarrow H'$ and
$\gamma:K\Rightarrow K'$ induces a natural transformation :
$$
\sCart_\cB(\beta,\gamma):\sCart_\cB(H,K)\Rightarrow
\sCart_\cB(H',K')
\qquad
G\mapsto\gamma*G*\beta
$$
and again, if $\beta=\one_H$ (resp. $\gamma=\one_K$) we also
denote this natural transformation by $\sCart_\cB(H,\gamma)$
(resp. $\sCart_\cB(\alpha,K)$).

(iii)\ \
Let $\sU,\sV$ be two universes, with $\sU\subset\sV$; we say
that the fibration $\phi$ has {\em essentially $\sU$-small
fibres} if $\phi^{-1}B$ is an essentially $\sU$-small category
for every $B\in\Ob(\cB)$. The $\sV$-small fibrations over $\cB$
with essentially $\sU$-small fibres form a $2$-category
$$
(\sU,\sV)\tdu\Fib(\cB)
$$
with $\Hom$-category $\sCart_\cB(\cC,\cC')$, for every pair
of fibrations $\cC\to\cB\leftarrow\cC'$, and with composition
functors given by (ii). When there is no danger of ambiguities,
we shall often write $\sU\tdu\Fib(\cB)$, or even just $\Fib(\cB)$
for this $2$-category. Notice that if $\cB$ is essentially
$\sU$-small, then $(\sU,\sV)\tdu\Fib(\cB)$ has essentially
$\sU$-small $\Hom$-categories.
\end{definition}

\begin{remark}\label{rem_added-little-extra}
Let $\cB$ and $\cB'$ be two categories and $\psi:\cB'\to\cB$
any functor.

(i)\ \
For every category $\cA$ and every functor $\phi:\cA\to\cB$,
the cartesian morphisms of the projection
$\phi':\cA\times_{(\phi,\psi)}\cB'\to\cB'$ are the pairs $(f,f')$
where $f$ is a cartesian morphism of $\cA$, $f'$ is a morphism
of $\cB'$, and $\phi f=\psi f'$ (see example \ref{ex_cat-cats}(i)).
Moreover, the projection $\cA\times_{(\phi,\psi)}\cB'\to\cA$ restricts
to a natural isomorphism of fibre categories
$$
\phi'^{-1}B'\isom\phi^{-1}(\psi B')
\qquad
\text{for every $B'\in\Ob(\cB')$}.
$$
If $\phi$ is a fibration, the same then holds for $\phi'$. Therefore,
for every pair of universes $\sU\subset\sV$ such that $\cB$ and
$\cB'$ are $\sV$-small, we get a well defined strict pseudo-functor
$$
(\sU,\sV)\tdu\Fib(\psi)^*:
(\sU,\sV)\tdu\Fib(\cB)\to(\sU,\sV)\tdu\Fib(\cB')
\qquad
(\cA\xrightarrow{\phi}\cB)\mapsto
(\cA\times_{(\phi,\psi)}\cB'\xrightarrow{\phi'}\cB')
$$
which assigns to every $\cB$-cartesian functor $F:\cA_1\to\cA_2$
the $\cB'$-cartesian functor
$F\times_\cB\cB':\cA_1\times_\cB\cB'\to\cA_2\times_\cB\cB'$, and
to any natural $\cB$-transformation $\alpha:F_1\Rightarrow F_2$
between $\cB$-cartesian functors $F_1,F_2:\cA_1\to\cA_2$ the
induced $\cB'$-transformation
$\alpha\times_\cB\cB':F_1\times_\cB\cB'\Rightarrow F_2\times_\cB\cB'$.
As usual, we often write $\sU\tdu\Fib(\psi)^*$, or just $\Fib(\psi)^*$
instead of $(\sU,\sV)\tdu\Fib(\psi)^*$.

(ii)\ \
Suppose that $\psi$ is a fibration with essentially $\sU$-small
fibres; then for every fibration $F':\cA'\to\cB'$ with essentially
$\sU$-small fibres, the composition $\psi\circ F':\cA'\to\cB$ is
also a fibration with essentially $\sU$-small fibres, and in this
case we get therefore also a well defined strict pseudo-functor
$$
(\sU,\sV)\tdu\Fib(\cB')\to(\sU,\sV)\tdu\Fib(\cB)
\qquad
(F':\cA'\to\cB')\mapsto(\psi\circ F:\cA'\to\cB).
$$

(iii)\ \
For $i=1,2$, let $\cA_i\to\cB$ be two fibrations; combining
(i) and (ii) we see that the induced functor
$\cA_1\times_\cB\cA_2\to\cB$ is also a fibration.

(iv)\ \
In the situation of (i), it is easily seen that
$\sCart_{\cB'}(\cB',\cA\times_{(\phi,\psi)}\cB')$ is naturally
identified  with the full subcategory of $\sCart_\cB(\cB',\cA)$
whose objects are the $\cB$-functors $\cB'\to\cA$ that send
{\em every} morphism of $\cB'$ to a $\cB$-cartesian morphism
of $\cA$ : the detailed verification shall be left to the reader.

(v)\ \
Suppose that $F:\cA\to\cA'$ is a {\em $\cB$-equivalence},
{\em i.e.} an equivalence in the $2$-category $\bCat/\cB$
(in the sense of definition \ref{def_adjoint-1-cells}(iii)).
Then it is easily seen that $F$ is $\cB$-cartesian.

(vi)\ \
In the situation of (i), let $\psi':\cB''\to\cB'$ be another
functor between $\sV$-small categories; then we have a natural
isomorphism of categories :
$$
(\sU,\sV)\tdu\Fib(\psi')^*\circ(\sU,\sV)\tdu\Fib(\psi)^*(\cA)\isom
(\sU,\sV)\tdu\Fib(\psi\circ\psi')^*(\cA)
$$
that assigns to every object $(A,B',B'')$ of
$(\cA\times_{(\phi,\psi)}\cB')\times_{(\phi',\psi')}\cB''$ the
object $(A,B'')$ of $\cA\times_{(\phi,\psi\circ\psi')}\cB''$,
and it is likewise defined on morphisms. Clearly this
isomorphism is strictly pseudo-natural with respect to $\cA$,
so we get a strict pseudo-natural isomorphism of strict
pseudo-functors :
$$
\gamma^{(\sU,\sV)-\Fib}_{\psi',\psi}:
(\sU,\sV)\tdu\Fib(\psi')^*\circ(\sU,\sV)\tdu\Fib(\psi)^*\isom
(\sU,\sV)\tdu\Fib(\psi\circ\psi')^*.
$$
\end{remark}

\sset\subsubsection{Cleavages}\label{subsec_choice-of-pseudo}
Let $\phi:\cA\to\cB$ be a fibration, and for every object
$(A,g:B\to FA)$ of $\cB/\phi\cA$ let us choose an inverse
image of $A$ over $g$ :
$$
g_A:g^*A\to A.
$$
Then we claim that the rule $(A,g)\mapsto g_A$ extends
uniquely to a functor
$$
\blambda:\cB/\phi\cA\to\sMorph(\cA)
$$
that makes commute the resulting diagram (notation of
\eqref{subsec_Morph-cat}) :
\set\begin{equation}\label{eq_elevator}
{\diagram \cB/\phi\cA \ar[rr]^-\blambda
\ar[rd]_\sS & & \sMorph(\cA) \ar[ld]^{\sMorph(\phi)} \\
& \sMorph(\cB).
\enddiagram}
\end{equation}
Indeed, since $g_A$ is $\phi$-cartesian, for every morphism
$(h,k):(A',B'\xrightarrow{g'}\phi A')\to(A,B\xrightarrow{g}\phi A)$
of $\cB/\phi\cA$ there exists a unique morphism $f:g'^*A'\to g^*A$
that makes commute the diagram
\set\begin{equation}\label{eq_square-elevator}
{\diagram g'^*A' \ar[r]^-{g'_{A'}} \ar[d]_f & A' \ar[d]^h \\
g^*A \ar[r]^-{g_A} & A
\enddiagram}
\qquad\text{and such that}\qquad
\phi(f)=k
\end{equation}
and this square diagram can be regarded as a morphism
$(f,h):g'_{A'}\to g_A$ in $\sMorph(\cA)$ such that
$\sMorph(\phi)(f,h)=(k,\phi(h))=\sS(h,k)$. The uniqueness
of $f$ easily implies that the rule $(h,k)\mapsto(f,h)$
yields a well defined functor as sought. We call a
{\em cleavage} for $\phi$ any such functor ({\em ``clivage''}
in french); hence, the cleavages for $\phi$ are characterized
as the functors $\blambda$ that make commute
\eqref{eq_elevator} and that map every object of
$\cB/\phi\cA$ to a $\phi$-cartesian morphism.

Notice that we may always choose a cleavage such that
$\blambda(A,\one_{FA})=\one_A$ for every $A\in\Ob(\cA)$.
We call {\em unital} a cleavage fulfilling this condition.

\begin{example}\label{ex_Morph-is-fibred}
Let $\phi:\cA\to\cB$ be any fibration. Then
$\sMorph(\phi):\sMorph(\cA)\to\sMorph(\cB)$ is a
fibration as well (notation of \eqref{subsec_Morph-cat}).
Indeed, pick a cleavage $\blambda$ for $\phi$; consider
any object $(A'\xrightarrow{h}A)$ of $\sMorph(\cA)$
and a morphism $(g',g):(B'\xrightarrow{k}B)\to
(\phi A'\xrightarrow{\phi(h)}\phi A)$ of $\sMorph(\cB)$.
By definition, $g:B\to\phi A$ and $g':B'\to\phi A'$ are
morphisms of $\cB$ such that $\phi(h)\circ g'=g\circ k$,
hence $(h,k):(A',B'\xrightarrow{g'}\phi A')\to
(A,B\xrightarrow{g}\phi A)$ is a morphism of $\cB/\phi\cA$,
and $(f,h):=\blambda(h,k)$ is a commutative square
\eqref{eq_square-elevator}. With this notation, it is
easily seen that $(g'^*A'\xrightarrow{f}g^*A)$ is an inverse
image of $(A'\xrightarrow{h}A)$ over $(g',g)$ : the details
are left to the reader.
\end{example}

\sset\subsubsection{The pseudo-functor associated with
a cleavage}\label{subsec_associate-c}
Let now $\blambda$ be any cleavage for $\phi$; to any
morphism $g:B'\to B$ in $\cB$ we attach a functor
$$
\sc_g:\phi^{-1}B\to\phi^{-1}B'
$$
as follows. First, we consider the functor
$$
(-,g):\phi^{-1}B\to B'/\phi\cA
\qquad
A\mapsto(A,g)
\qquad
(A_1\xrightarrow{\ h\ }A_2)\mapsto
((A_1,g)\xrightarrow{\ (h,\one_{B'})\ }(A_2,g)).
$$
Let also $\ss:\sMorph(\cA)\to\cA$ be the source functor
(see \eqref{subsec_Morph-cat}), and notice that the
composition $\ss\circ\blambda\circ(-,g)$ factors through
$\phi^{-1}B'$, and yields therefore the sought functor
$\sc_g$. Moreover, the commutativity of the diagram
\eqref{eq_square-elevator} also means that the rule :
$A\mapsto\blambda(A,g)$ defines a natural transformation
\set\begin{equation}\label{eq_associated-pseudo-cocone}
g_\bullet:\iota_{B'}\circ\sc_g\Rightarrow\iota_B
\end{equation}
(notation of definition \ref{def_pull-bekko}(i)). Notice
especially that every $B\in\Ob(\cB)$ yields a natural
isomorphism of functors
$$
\delta_B:\one_{\phi^{-1}B}\Rightarrow\sc_{\one_B}
\qquad
A\mapsto((\one_B)^{-1}_A:A\isom(\one_B)^*A).
$$
Furthermore, let $B''\xrightarrow{h} B'\xrightarrow{g} B$
be two morphisms in $\cB$, so we get the functors
$$
\sc_g:\phi^{-1}B\to\phi^{-1}B' \qquad
\sc_h:\phi^{-1}B'\to\phi^{-1}B'' \qquad
\sc_{g\circ h}:\phi^{-1}B\to\phi^{-1}B''
$$
as well as the natural transformations :
$$
g_\bullet:\iota_{B'}\circ\sc_g\Rightarrow\iota_B \qquad
h_\bullet:\iota_{B''}\circ\sc_h\Rightarrow\iota_{B'} \qquad
(g\circ h)_\bullet:\iota_{B''}\circ\sc_{g\circ h}\Rightarrow\iota_B.
$$
By inspecting the constructions, one easily finds a unique
natural isomorphism :
$$
\gamma_{(h,g)}:\sc_h\circ\sc_g\Rightarrow\sc_{g\circ h}
$$
which fits into a commutative diagram :
$$
\xymatrix{
\iota_{B''}\circ\sc_h\circ\sc_g
\ar@{=>}[rr]^-{h_\bullet*\sc_g} \ar@{=>}[d]_{\iota_{B''}*\gamma_{(h,g)}}
& & \iota_{B'}\circ\sc_g \ar@{=>}[d]^{g_\bullet} \\
\iota_{B''}\circ\sc_{g\circ h} \ar@{=>}[rr]^-{(g\circ h)_\bullet} & &
\iota_B.
}$$
Moreover, if $k:B'''\to B''$ is a third morphism of $\cB$,
we can compute :
$$
\begin{aligned}
(g\circ h\circ k)_A\circ\gamma_{(k,g\circ h),A}\circ\sc_k(\gamma_{(h,g),A})
=&\,(g\circ h)_A\circ k_{\sc_{g\circ h}A}\circ\sc_k(\gamma_{(h,g),A}) \\
=&\,(g\circ h)_A\circ\gamma_{(h,g),A}\circ k_{\sc_h\sc_gA} \\
=&\,g_A\circ h_{\sc_gA}\circ k_{\sc_h\sc_gA} \\
=&\,g_A\circ(h\circ k)_{\sc_gA}\circ\gamma_{(k,h),\sc_gA} \\
=&\,(g\circ h\circ k)_A\circ\gamma_{(h\circ k,g),A}\circ\gamma_{(k,h),\sc_gA}
\end{aligned}
$$
for every $A\in\Ob(\cA)$, and since $(g\circ h\circ k)_A$
is $\cB$-cartesian, we deduce :
\set\begin{equation}\label{eq_tostaki}
\gamma_{(k,g\circ h),A}\circ\sc_k(\gamma_{(h,g),A})=
\gamma_{(h\circ k,g),A}\circ\gamma_{(k,h),\sc_gA}.
\end{equation}
Likewise, for every morphism $f:B'\to B$ of $\cB$ and
every $A\in\Ob(\phi^{-1}B)$ we have the commutative diagram :
$$
\xymatrix{ f^*\one^*_BA \ar[rr]^-{\sc_f(\one_B)_A}
\ar[d]_{f_{\one^*_BA}} & & f^*A \ar[d]^{f_A} & &
\one^*_{B'}f^*A \ar[ll]_-{(\one_{B'})_{f^*A}}
\ar[d]^{\gamma_{(\one_{B'},f),A}} \\
\one^*_BA \ar[rr]^-{(\one_B)_A} & & A & &
f^*A \ar[ll]_-{f_A}
}$$
which implies :
\set\begin{equation}\label{eq_nie-en-bloc}
\gamma_{(f,\one_B),A}=\sc_f((\one_B)_A)
\qquad
\gamma_{(\one_{B'},f),A}=(\one_{B'})_{\sc_fA}.
\end{equation}
Suppose now that $\phi$ has $\sV$-small fibres for some
universe $\sV$. Then the identities \eqref{eq_tostaki} and
\eqref{eq_nie-en-bloc} mean that the rule which assigns to
each $B\in\Ob(\cB)$ the small category $\phi^{-1}B$ and to
each morphism $g$ in $\cB$ the functor $\sc_g$ defines a
pseudo-functor
$$
\sc:\cB^o\to\sV\tdu\bCat
$$
whose coherence constraint is the system of natural
isomorphisms $(\delta_\bullet,\gamma_\bullet)$. Moreover,
the system of inclusion functors $(\iota_B~|~B\in\Ob(\cB))$
amounts to a pseudo-cocone
$$
\iota:\sc\Rightarrow\sF_\cA
$$
whose coherence constraint is given by the system of
natural transformations \eqref{eq_associated-pseudo-cocone}. 
Notice that if $\blambda$ is a unital cleavage, $\sc$
will be a unital pseudo-functor. (Here we view $\cB^o$ as a
$2$-category, as explained in example \ref{ex_2-cats}(i);
also, the $2$-category structure on $\sV\tdu\bCat$ is the
one of remark \ref{rem_equiv-2-cat}(ii)). We call $\sc$ and
$\iota$ respectively {\em the pseudo-functor and the
pseudo-cocone associated with the cleavage $\blambda$}.

\begin{example}\label{ex_Fubini-fibred}
As an application, we may generalize example
\ref{ex_lim_interchange} as follows. Let $I,J$ be two small
categories, $\phi:I\to J$ a fibration, fix a cleavage
$\blambda$ for $\phi$, and denote by $\sc$ the associated
pseudo-functor.

(i)\ \
First, we claim that for every $j\in\Ob(J)$ the functor
$(-,\one_j):\phi^{-1}j\to j/\phi I$ is final (notation of
\eqref{subsec_associate-c}). Indeed, for any
$(i_0,f:j\to\phi i_0)\in\Ob(j/\phi I)$ we have the
cartesian morphism $\blambda(i_0,f):\sc_f(i_0)\to i_0$,
and it is easily seen that the category
$(-,\one_j)(\phi^{-1}j)/(i_o,f)$ is isomorphic to
$(\phi^{-1}j)/\sc_f(i_0)$, which is obviously connected,
whence the claim.

(ii)\ \
Notice that the composition
$\st_j\circ(-,\one_j):\phi^{-1}j\to I$ equals the inclusion
functor $\iota_j$ (here $\st_j:j/\phi I\to I$ is the target
functor). Moreover, for every morphism $f:j\to j'$ in $J$,
the natural transformation
$f_\bullet:\iota_j\circ\sc_f\Rightarrow\iota_{j'}$ associated
with $\blambda$ induces a morphism
$$
\omega^f_F:\lim_{\phi^{-1}j}F\circ\iota_j
\xrightarrow{\ \ \underset{\sc_f}{\mathrm{lim}}\,\one_{F\circ\iota_j}\ \ }
\lim_{\phi^{-1}{j'}}F\circ\iota_j\circ\sc_f
\xrightarrow{\ \ \underset{\phi^{-1}{j'}}{\mathrm{lim}}\,F*f_\bullet\ \ }
\lim_{\phi^{-1}{j'}}F\circ\iota_{j'}
$$
for every functor $F:I\to\cC$, fitting into a commutative diagram
$$
\xymatrix{
\int_\phi^\wedge F(j) \ar[rr]^-{\int^\wedge_\phi F(f)} \ar[d]
& & \int^\wedge_\phi F(j') \ar[d] \\
\displaystyle\lim_{\phi^{-1}j}F\circ\iota_j
\ar[rr]^-{\omega^f_F} & &
\displaystyle\lim_{\phi^{-1}{j'}}F\circ\iota_{j'}
}$$
whose vertical arrows are the natural isomorphisms
provided by (i) and remark \ref{rem_fun-cofinal}(ii).
In other words, $\int^\wedge_\phi$ is naturally isomorphic
to the functor that assigns to every such $F$ the functor
$$
\int^\wedge_\blambda F:J\to\cC^\wedge
\qquad
j\mapsto\lim_{\phi^{-1}j}F\circ\iota_j
\qquad
(f:j\to j')\mapsto\omega^f_F.
$$

(iii)\ \
Furthermore, in this situation, proposition \ref{prop_Fubini}
says that there is a natural isomorphism
$$
\lim_IF\isom\Lim_J\int^\wedge_\blambda F
\qquad
\text{in $\cC^\wedge$}
$$
for every functor $F:I\to\cC$. If $\cC$ is complete, the
functor $\int^\wedge_\blambda$ is isomorphic, by remark
\ref{rem_interchange}, to the composition of $\bFun(J,h_\cC)$
and a functor well defined up to isomorphism
$$
\int_\blambda:\bFun(I,\cC)\to\bFun(J,\cC)
\qquad
F\mapsto(j\mapsto\Lim_{\phi^{-1}j}F\circ\iota_j)
$$
and we can restate the foregoing ``Fubini'' isomorphism in
terms of this latter functor.

(iv)\ \
Lastly, let $G:I^o\to\cC$ be any functor. It follows from
(ii) that the functor $\int_\wedge^{\phi^o}G$ is naturally
isomorphic to the functor
$$
\int^\blambda_\wedge G:J^o\to\cC^{o\wedge o}
\qquad
j^o\mapsto\colim_{(\phi^{-1}j)^o}G\circ\iota_j^o
\qquad
(f^o:j'^o\to j^o)\mapsto\omega^f_{G^o}
$$
and we have a natural isomorphism
$$
\colim_{I^o}G\isom\colim_{J^o}\int^\blambda_\wedge G
\qquad
\text{in $\cC^{o\wedge o}$}.
$$
Again, if $\cC$ is cocomplete, the functor $\int^\blambda_\wedge$
is the composition of $\bFun(J^o,h^o_{\cC^o})$ and a functor
$$
\int^\blambda:\bFun(I^o,\cC)\to\bFun(J^o,\cC)
\qquad
G\mapsto(j^o\mapsto\Colim_{(\phi^{-1}j)^o}G\circ\iota^o_j)
$$
and we may state the foregoing ``Fubini'' isomorphism for colimits
in terms of this latter functor.
\end{example}

\sset\subsubsection{Fibration associated with a presheaf}
\label{subsec_fibred-cats-II}
Let $\cB$ be a category, $F$ a presheaf on $\cB$. As in
\eqref{subsec_category-of-elements}, we let $\cFib(F)$ be
the category of elements of $F$, and we notice that the
source functor
$$
\ss_F:\cFib(F)\to\cB
$$
is a fibration. For every $X\in\Ob(\cB)$, the fibre
$\ss_F^{-1}(X)$ is (naturally isomorphic to) the discrete
category $FX$ (see example \ref{ex_universe}(ii)). Notice
also that every morphism in $\cFib(F)$ is cartesian, and for
every $\cB$-category $\cC$, the category
$\sCart_\cB(\cC,\cFib(F))$ is discrete.

For every pair of presheaves $F,G$ on $\cB$, we have a natural
bijection:
$$
\Hom_{\cC^\wedge}(F,G)\isom\Ob(\sCart_\cB(\cFib(F),\cFib(G)))
$$
(which we can view as an isomorphism of discrete categories).
Namely, to a morphism $\psi:F\to G$ one assigns the functor
$$
\cFib(\psi):\cFib(F)\to\cFib(G)
\qquad
(X,s)\mapsto(X,\psi_X(s))
\quad
\text{for every $(X,s)\in\Ob(\cFib(F))$}.
$$

\begin{example}\label{ex_fibred-cats-II}
(i)\ \
For instance, if $\cC$ has small $\Hom$-sets, then one checks
easily that there is a natural isomorphism of $\cC$-fibrations
$$
\cFib(h_X)\isom\cC/X
\qquad
\text{for every $X\in\Ob(\cC)$}
$$
where $h_X$ is the presheaf represented by $X$ (see
\eqref{subsec_yoneda}), and $\cC/X$ is regarded as a
$\cC$-fibration as in example \ref{ex_fibred-cats}(i).
Moreover, every morphism $f:X\to Y$ of $\cC$ induces
a morphism of presheaves $h_f:h_Y\to h_X$, and the
foregoing isomorphism identifies $\cFib(h_f)$ with the
functor $f_*:\cC/X\to\cC/Y$ of \eqref{eq_push-for}.

(ii)\ \
In the situation of (i), let $u:\cC'\to\cC$ be any functor;
by a direct inspection we get a natural isomorphism of
categories
$$
\Fib(u)^*(\cC/X)\isom u\cC'/X
$$
(notation of remark \ref{rem_added-little-extra}(i)) which
identifies the fibration $\Fib(u)^*(\cC/X)\to\cC'$ with the
source functor $\ss_X:u\cC'/X\to\cC'$ of
\eqref{subsec_fibreovercat}. Likewise, for every morphism
$f:X\to Y$ in $\cC$, the cartesian functor
$\Fib(u)^*(f_*):\Fib(u)^*(\cC/X)\to\Fib(u)^*(\cC/Y)$ is
identified with the functor $u\cC'/f:u\cC'/X\to u\cC'/Y$
as in \eqref{subsec_fibreovercat}.

(iii)\ \
As a special case, let $h_\cC:\cC\to\cC^\wedge$ be the Yoneda
embedding, and $F$ any presheaf on $\cC$. As noted in
\eqref{subsec_category-of-elements}, we have a natural
isomorphism of categories $h_\cC\cC/F\isom\cFib(F)$, so (i)
generalizes to a natural isomorphism of $\cC$-fibrations :
$$
\cFib(F)\isom\sV\tdu\Fib(h_\cC)^*(\cC^\wedge/F)
$$
for every universe $\sV$ such that $\cC^\wedge$ has $\sV$-small
$\Hom$-sets.
\end{example}

\sset\subsubsection{}\label{subsec_conversely}
Conversely, let $\phi:\cA\to\cB$ be a fibration {\em with
discrete fibres}, {\em i.e.} such that $\phi^{-1}B$ is a
discrete small category for every $B\in\Ob(\cB)$. Then it is
easily seen that for every morphism $g:B'\to B$ of $\cB$,
every $A\in\phi^{-1}B$ admits a unique inverse image over
$g$. It follows that $\phi$ admits a unique cleavage.
Moreover the  associated pseudo-functor $\sc:\cB^o\to\bCat$
is strict, hence it may be regarded as a functor with values
in the subcategory $\Set$ of $\bCat$; {\em i.e.} $\sc$ is
a presheaf on $\cB$, and a simple inspection shows that the
fibration $\cFib(\sc)$ is naturally isomorphic to $\phi$.
Summing up, we have obtained a fully faithful functor
$$
\cFib:\cB^\wedge\to\Fib(\cB)
$$
whose essential image is the full subcategory of $\Fib(\cB)$
whose objects are the fibrations with discrete fibres. We shall
show next how to extend this equivalence to the whole of $\Fib(\cB)$.

\sset\subsubsection{Fibration associated with a pseudo-functor}
\label{subsec_fib-from-pseudo}
Let $\cB$ be any category, and $\sc:\cB^o\to\bCat$ any
pseudo-functor, with coherence constraint $(\delta^\sc,\gamma^\sc)$.
We attach to $\sc$ the $\cB$-category
$$
\pi^\sc:\cFib(\sc)\to\cB
$$
such that
$\Ob(\cFib(\sc)):=\{(B,X)~|~B\in\Ob(\cB),\ X\in\Ob(\sc_B)\}$,
and where
$$
\Hom_{\cFib(\sc)}((B,X),(B',Y)):=
\{(\phi,f)~|~\phi\in\Hom_\cB(B,B'),\ f\in\Hom_{\sc_B}(X,\sc_\phi Y)\}
$$
for every two objects $(B,X),(B',Y)$. The composition
of two morphisms $f:X\to\sc_\phi Y$ and $g:Y\to\sc_\psi Z$
is the pair $(\psi\circ\phi,t)$, where $t$ is the composition
$$
X\xrightarrow{\ f\ }\sc_\phi Y\xrightarrow{\ \sc_\phi(g)\ }
\sc_\phi\sc_\psi Z\xrightarrow{\ \gamma^\sc_{(\phi,\psi),Z}\ }
\sc_{\psi\circ\phi}Z.
$$
The unit axiom for $\delta^\sc$ implies easily that for every
object $(B,X)$ the morphism
$$
(\one_B,\delta^\sc_{B,X}:X\to\sc_{\one_B}X)
$$
is neutral for left and right composition with any other
morphism of $\cFib(\sc)$. Let us show the associativity : if
$h:Z\to\sc_\lambda W$ is a third morphism, we need to verify
the identity
$$
\gamma^\sc_{(\psi\circ\phi,\lambda),W}\circ\sc_{\psi\circ\phi}(h)\circ
\gamma^\sc_{(\phi,\psi),Z}\circ\sc_\phi(g)\circ f=
\gamma^\sc_{(\phi,\lambda\circ\psi),W}\circ\sc_\phi(\gamma^\sc_{(\psi,\lambda),W})
\circ\sc_\phi\sc_\psi(h)\circ\sc_\phi(g)\circ f.
$$
But we have $\sc_{\psi\circ\phi}(h)\circ\gamma^\sc_{(\phi,\psi),Z}=
\gamma^\sc_{(\phi,\psi),\sc_\lambda W}\circ\sc_\phi\sc_\psi(h)$
by the naturality of $\gamma^\sc_{(\phi,\psi)}$, hence we are
reduced to checking that
$$
\gamma^\sc_{(\psi\circ\phi,\lambda),W}\circ\gamma^\sc_{(\phi,\psi),\sc_\lambda W}=
\gamma^\sc_{(\phi,\lambda\circ\psi),W}\circ\sc_\phi(\gamma^\sc_{(\psi,\lambda),W})
$$
which follows from the composition axioms for $\gamma^\sc$.
The functor $\pi^\sc$ is given by the rules : $(B,X)\mapsto B$
and $(\phi,f)\mapsto\phi$ for every object $(B,X)$ and
every morphism $(\phi,f)$ of $\cFib(\sc)$. 

\sset\subsubsection{}\label{subsec_sd}
Let $\sd:\cB^o\to\bCat$ be another pseudo-functor, with
coherence constraint $(\delta^\sd,\gamma^\sd)$, and
$$
\omega:\sc\Rightarrow\sd
$$
any pseudo-natural transformation, with coherence constraint
$\tau^\omega$. We define a $\cB$-functor
$$
\cFib(\omega):\cFib(\sc)\to\cFib(\sd)
$$
by assigning to every object $(B,X)$ of $\cFib(\sc)$ the
object $(B,\omega_BX)$ of $\cFib(\sd)$, and to every morphism
$(\phi,f):(B,X)\to(B',Y)$ the morphism $(\phi,t)$, where $t$
is the composition
$$
\omega_BX\xrightarrow{\ \omega_Bf\ }\omega_B\circ\sc_\phi Y
\xrightarrow{\ (\tau^\omega_{\phi,Y})^{-1}\ }\sd_\phi\circ\omega_{B'}Y.
$$
To check that these rules do yield a functor, consider any
other morphism $(\psi,g):(B',Y)\to(B'',Z)$ of $\cFib(\sc)$;
we need to show that
$$
\tau^{\omega\ -1}_{\psi\circ\phi,Z}\circ
\omega_B(\gamma^\sc_{(\phi,\psi),Z}\circ\sc_\phi(g)\circ f)=
\gamma^\sd_{(\phi,\psi),\omega_{B''}Z}\circ\sd_\phi(\tau^{\omega\ -1}_{\psi,Z}
\circ\omega_{B'}(g))\circ\tau^{\omega\ -1}_{\phi,Y}\circ\omega_B(f).
$$
However, the naturality of $\tau^\omega_\phi$ implies that
$\tau^{\omega\ -1}_{\phi,\sc_\psi Z}\circ\omega_B(\sc_\phi g)=
\sd_\phi(\omega_{B'}g)\circ\tau^{\omega\ -1}_{\phi,Y}$,
so we are reduced to checking that :
$$
\tau^{\omega\ -1}_{\psi\circ\phi,Z}\circ\omega_B(\gamma^\sc_{(\phi,\psi),Z})=
\gamma^\sd_{(\phi,\psi),\omega_{B''}Z}\circ\sd_\phi(\tau^{\omega\ -1}_{\psi,Z})
\circ\tau^{\omega\ -1}_{\phi,\sc_\psi Z}.
$$
But the latter follows directly from the coherence axioms
for $\tau^\omega$. Likewise, the assertion that $\cFib(\omega)$
respects identity morphisms comes down to the equalities :
$$
\tau^{\omega\ -1}_{\one_B,X}\circ\omega_B(\delta^\sc_{B,X})=
\delta^\sd_{B,\omega_BX}
$$
which again follows from the coherence axioms for $\tau^\omega$.
Moreover, if $\se:\cB^o\to\bCat$ is a third pseudo-functor, and
$$
\mu:\sd\Rightarrow\se
$$
a second pseudo-natural transformation, we have
$$
\cFib(\mu)\circ\cFib(\omega)=\cFib(\mu\odot\omega).
$$
Indeed, if $\tau^\mu$ and $\tau^{\mu\odot\omega}$ are the coherence
constraints for $\mu$ and $\mu\odot\omega$, the assertion
amounts to :
$$
\tau^{\mu\ -1}_{\phi,\omega_{B'}Y}\circ\mu_B(\tau^\omega_{\phi,Y})^{-1}
\circ\mu_B(\omega_Bf)=\tau^{\mu\odot\omega\ -1}_{B,Y}\circ\mu_B(\omega_Bf)
$$
for every morphism $(\phi,f):(B,X)\to(B',Y)$ in $\cFib(\sc)$,
which is clear from the definition of $\mu\odot\omega$.

\begin{lemma}\label{lem_distinguished-cleavage}
With the notation of \eqref{subsec_sd}, the following holds :
\begin{enumerate}
\item
A morphism $(\phi,f):(B,X)\to(B',Y)$ of $\cFib(\sc)$ is
$\pi^\sc$-cartesian if and only if $f:X\to\sc_\phi Y$ is
an isomorphism in $\sc_B$.
\item
The functor $\pi^\sc$ is a fibration, and is endowed with
a distinguished cleavage
$$
\blambda^*:\cB/\pi^\sc\cFib(\sc)\to\sMorph(\cFib(\sc))
\qquad
(X,\phi)\mapsto(\phi,\one_{\sc_\phi X}).
$$
\item
The functor $\cFib(\omega)$ is $\cB$-cartesian.
\end{enumerate}
\end{lemma}
\begin{proof}(i): In view of proposition
\ref{prop_towards-2-yoneda} and the discussion of
\eqref{subsec_sd}, we may assume that $\blambda$ is
unital (details left to the reader). Now, suppose
first that $(\phi,f)$ is $\cB$-cartesian; then there
exists a unique morphism $(\one_B,g):(B,\sc_\phi Y)\to(B,X)$
such that $(\phi,f)\circ(\one_\B,g)=(\phi,\one_{\sc_\phi Y})$.
Since $\sc$ is unital, we easily deduce that
$f\circ g=\one_{\sc_\phi Y}$. Moreover,
$(\one_B,g\circ f):(B,X)\to(B,X)$ is a morphism of
$\cFib(\sc)$ such that
$(\phi,f)\circ(\one_B,g\circ f)=(\phi,f)=
(\phi,f)\circ(\one_B,\one_X)$, whence $g\circ f=\one_X$.
Conversely, suppose that $f$ is an isomorphism, and let
$\psi:C\to B$ any morphism in $\cB$, $Z\in\Ob(\sc_C)$
any object, and $(\phi\circ\psi,h):(C,Z)\to(B',Y)$
any morphism of $\cFib(\sc)$; we let $k:Z\to\sc_\psi X$
be the composition
$$
Z\xrightarrow{\ h\ }\sc_{\phi\circ\psi}Y
\xrightarrow{\ \gamma^{\sc\ -1}_{(\psi,\phi),Y}\ }
\sc_\psi\sc_\phi Y\xrightarrow{\ \sc_\phi(f^{-1})\ }\sc_\psi X.
$$
Then $(\psi,k):(C,Z)\to(B,X)$ is the unique morphism of
$\cFib(\sc)$ such that $(\phi,f)\circ(\psi,k)=(\phi\circ\psi,h)$;
this shows that $(\phi,f)$ is $\cB$-cartesian.

(ii) and (iii) are immediate consequences of (i).
\end{proof}

\sset\subsubsection{}\label{subsec_from-mod-to-trans}
Lastly, consider pseudo-natural transformations
$$
\omega,\omega':\sc\Rightarrow\sd
\qquad\text{and a modification}\qquad
\Xi:\omega\leadsto\omega'.
$$
We attach to $\Xi$ the natural transformation of functors
$$
\cFib(\Xi):\cFib(\omega)\Rightarrow\cFib(\omega')
$$
assigning to every object $(B,X)$ of $\cFib(\sc)$ the
morphism $(\one_B,t)$ of $\cFib(\sd)$, where $t$ is the
composition
$$
\omega_BX\xrightarrow{\ \Xi_{B,X}\ }\omega'_BX
\xrightarrow{\ \delta^\sd_{B,\omega'_BX}\ }\sd_{\one_B}(\omega'_BX).
$$
In order to verify the naturality of $\cFib(\Xi)$, it
suffices to check the commutativity of the diagram:
$$
\xymatrix{ \omega_BX \ar[rrr]^-{\Xi_{B,X}} \ar[d]_{\omega_Bf}
& & & \omega'_BX \ar[rrr]^-{\delta^\sd_{B,\omega'_BX}}
\ar[d]_{\omega'_Bf} & & &
\sd_{\one_B}(\omega'_BX) \ar[d]^{\sd_{\one_B}(\omega'_Bf)} \\
\omega_B\sc_\phi Y \ar[rrr]^-{\Xi_{B,\sc_\phi Y}}
\ar[d]_{\tau^{\omega\ -1}_{\phi,Y}} & & &
\omega'_B\sc_\phi Y \ar[rrr]^-{\delta^\sd_{B,\omega'_B\sc_\phi Y}}
\ar[d]_{\tau^{\omega'\ -1}_{\phi,Y}} & & & \sd_{\one_B}(\omega'_B\sc_\phi Y)
\ar[d]^{\gamma^\sd_{(\one_B,\phi),\omega'_{B'}Y}\circ\sd_{\one_B}(\tau^{\omega'\ -1}_{\phi,Y})} \\
\sd_\phi\omega_{B'}Y \ar[rrr]^-{\sd_\phi(\Xi_{B',Y})} & & &
\sd_\phi\omega'_{B'}Y
\ar[rrr]^-{\gamma^\sd_{(\phi,\one_B),\omega'_{B'}Y}\circ
\sd_\phi(\delta^\sd_{B',\omega'_{B'}Y})} & & & \sd_\phi\omega'_{B'}Y
}$$
for every morphism $(\phi,f):(B,X)\to(B',Y)$ of $\cFib(\sc)$.
However, the commutativity of the two top squares follows
from the naturality of $\Xi_B$ and $\delta^\sd_B$, and that
of the left bottom square translates the compatibility
condition for $\Xi$. Then, since
$$
\gamma^\sd_{(\one_B,\phi),\omega'_{B'}Y}=
\delta^{\sd\ -1}_{B',\sd_\phi\omega'_{B'}Y}
\qquad\text{and}\qquad
\gamma^\sd_{(\phi,\one_B),\omega'_{B'}Y}\circ
\sd_\phi(\delta^\sd_{B',\omega'_{B'}Y})=\one_{\sd_\phi\omega'_{B'}Y}
$$
we are reduced to checking the identity
$$
\sd_{\one_B}(\tau^{\omega'}_{\phi,Y})\circ\delta^\sd_{B',\sd_\phi\omega'_{B'}Y}
=\delta^\sd_{B,\omega'_B\sc_\phi Y}\circ\tau^{\omega'}_{\phi,Y}
$$
which follows from the naturality of $\delta^\sd_B$.
Furthermore, suppose we have a third pseudo-natural
transformation
$$
\omega'':\sc\Rightarrow\sd
\qquad\text{and a modification}\qquad
\Xi':\omega'\leadsto\omega''.
$$
Then we have as well
\set\begin{equation}\label{eq_compose-Xis}
\cFib(\Xi'\odot\Xi)=\cFib(\Xi')\odot\cFib(\Xi).
\end{equation}
Indeed, the assertion amounts to checking, for every
object $(B,X)$ of $\cFib(\sc)$, the identity
$$
\gamma^\sd_{(\one_B,\one_B),\omega''_BX}\circ
\sd_{\one_B}(\delta^\sd_{B,\omega''_BX})\circ
\sd_{\one_B}(\Xi'_{B,X})\circ\delta^\sd_{B,\omega'_BX}
\circ\Xi_{B,X}=
\delta^\sd_{B,\omega''_BX}\circ\Xi'_{B,X}\circ\Xi_{B,X}
$$
which is clear, since $\gamma^\sd_{(\one_B,\one_B),\omega''_BX}\circ
\sd_{\one_B}(\delta^\sd_{B,\omega''_BX})=\one_{\sd_{\one_B}\omega''_BX}$
by the unit axiom for $\delta^\sd$, and
$\sd_{\one_B}(\Xi'_{B,X})\circ\delta^\sd_{B,\omega'_BX}=
\delta^\sd_{B,\omega''_BX}\circ\Xi'_{B,X}$, by the naturality
of $\delta^\sd_B$. Likewise, suppose we have two pseudo-natural
transformations
$$
\mu,\mu':\sd\Rightarrow\se
\qquad\text{and a modification}\qquad
\Theta:\mu\leadsto\mu'.
$$
Then we have as well
$$
\cFib(\Theta)*\cFib(\Xi)=\cFib(\Theta*\Xi).
$$
For the proof, in light of \eqref{eq_compose-Xis} we may
assume that either $\Theta=\one_\mu$ or $\Xi=\one_\omega$.
In the first case, the assertion comes down to the identity
$$
(\tau^\mu_{\one_B,\omega'_BX})^{-1}\circ
\mu_B(\delta^\sd_{B,\omega'_BX}\circ\Xi_{B,X})=
\delta^\sd_{B,\mu_B\omega'_BX}\circ\mu_B(\Xi_{B,X})
$$
which is equivalent to :
$\mu_B(\delta^\sd_{B,\omega'_BX})=
\tau^\mu_{\one_B,\omega'_BX}\circ\delta^\sd_{B,\mu_B\omega'_BX}$, which
in turns follows from the coherence axioms for $\tau^\mu$.
The case where $\Xi=\one_\omega$ is clear by a simple
inspection. Summing up, we have obtained a strict
pseudo-functor
$$
\cFib_\cB:\sPsFun(\cB^o,\bCat)\to\Fib(\cB).
$$

\begin{remark}\label{rem_functors-from-lax}
(i)\ \
Notice that the discussion of \eqref{subsec_sd} requires
the invertibility of the coherence constraint $\tau^\omega$,
and therefore it does not apply to general lax-natural
transformations $\omega:\sc\Rightarrow\sd$. However,
with the obvious changes, it does apply to lax-natural
transformations $\omega:{}^o\sc\Rightarrow{}^o\sd$ (see
example \ref{ex_from-Cats-to-PsFuns}(ii)) : the details
shall be left to the reader. From lemma
\ref{lem_distinguished-cleavage} and its proof, we then
easily see that the resulting functor
$\cFib(\omega):\cFib(\sc)\to\cFib(\sd)$ will be cartesian
if and only if $\omega$ is pseudo-natural.

(ii)\ \
Likewise, the discussion of \eqref{subsec_from-mod-to-trans}
applies more generally to a modification
$\Xi:\omega'\leadsto\omega$ between any pair of lax-natural
transformations $\omega,\omega':{}^o\sc\Rightarrow{}^o\sd$;
namely, any such modification induces a natural transformation
$\cFib(\Xi):\cFib(\omega)\Rightarrow\cFib(\omega')$.
\end{remark}

\begin{theorem}\label{th_fundamental-fibrations}
For every category $\cB$, the pseudo-functor
$\cFib_\cB$ is a strong $2$-equivalence.
\end{theorem}
\begin{proof} Let $\phi:\cA\to\cB$ be any fibration with
small fibres; pick a unital cleavage $\blambda$ for $\phi$,
and let $\sc$ be the pseudo-functor associated with $\blambda$.

\begin{claim}\label{cl_exhibit-nat-isom}
There exists an isomorphism of $\cB$-categories
$$
\psi^\blambda:\cFib(\sc)\isom\cA
\qquad
(B,X)\mapsto X
\qquad
((B,X)\xrightarrow{\ (g,f)\ }(B',X'))\mapsto\blambda(X',g)\circ f.
$$
\end{claim}
\begin{pfclaim} Indeed, since $\sc$ is unital, it is clear
that $\psi^\blambda(\one_{(B,X)})=\one_{\psi(B,X)}$ for every
$(B,X)\in\Ob(\cFib(\sc))$. Next, let $(g,f):(B,X)\to(B',X')$
and $(g',f'):(B',X')\to(B'',X'')$ be two morphisms; then
$(g',f')\circ(g,f)=
(g'\circ g,\gamma^\sc_{(g,g'),X''}\circ\sc_g(f')\circ f)$,
and on the other hand, we may compute
$$
\begin{aligned}
\psi^\blambda(g',f')\circ\psi^\blambda(g,f)
=&\, \blambda(X'',g')\circ f'\circ\blambda(X',g)\circ f \\
=&\, \blambda(X'',g')\circ\blambda(\sc_{g'}X'',g)
     \circ\sc_g(f')\circ f \\
=&\, \blambda(X'',g'g)\circ\gamma^\sc_{(g,g'),X''}\circ
     \sc_g(f')\circ f \\
=&\, \psi^\blambda((g',f')\circ(g,f))
\end{aligned}
$$
as required. Clearly $\psi^\blambda$ is bijective on objects,
and since $\blambda(A,g)$ is a $\phi$-cartesian morphism for
every object $(A,g)$ of $\cB/\phi\cA$, it is easily seen
that $\psi^\blambda$ is fully faithful, hence it is an isomorphism
of $\cB$-categories.
\end{pfclaim}

In view of claim \ref{cl_exhibit-nat-isom}, it remains
only to check that for every pair of pseudo-functors
$\sc,\sd:\cB^o\to\bCat$ the induced functor
\set\begin{equation}\label{eq_sPsNat-to-cFib}
\sPsNat(\sc,\sd)\to\sCart_\cB(\cFib(\sc),\cFib(\sd))
\end{equation}
is an isomorphism of categories, and due to proposition
\ref{prop_towards-2-yoneda} and the discussion of
\eqref{subsec_sd}, we may assume that $\sc$ and $\sd$ are
unital. Thus, let $\omega,\omega':\sc\Rightarrow\sd$ be two
pseudo-natural transformations, with respective coherence
constraints $\tau^\omega$ and $\tau^{\omega'}$, and suppose
that $\cFib(\omega)=\cFib(\omega')$; then clearly
$\omega_BX=\omega'_BX$ for every object $(B,X)$ of
$\cFib(\sc)$. Likewise, in light of remark \ref{rem_unital}(ii)
we see that
$\tau^\omega_{\one_B,X}=\tau^{\omega'}_{\one_B,X}=\one_{\omega_BX}$
for every such $(B,X)$, whence $\omega_Bf=\omega'_Bf$
for every $B\in\Ob(\cB)$ and every morphism $f$ of $\sc_B$.
Then it is also clear that
$\tau^\omega_{\phi,Y}=\tau^{\omega'}_{\phi,Y}$ for every morphism
$\phi:B\to B'$ of $\cB$ and every $Y\in\Ob(\sc_{B'})$,
{\em i.e.} $\omega=\omega'$.

Next, let $F:\cFib(\sc)\to\cFib(\sd)$ be any $\cB$-cartesian
functor; for every object $(B,X)$ of $\cFib(\sc)$ we define
$\omega_BX$ as the unique object of $\sd_B$ such that
$F(B,X)=(B,\omega_BX)$. Likewise -- and since $\sd$ is unital --
for every morphism of $\cFib(\sc)$ of the form
$(\one_B,f):(B,X)\to(B,Y)$ we may define
$\omega_Bf:\omega_BX\to\omega_BY$ so that
$F(\one_B,f)=(\one_B,\omega_Bf)$. Since both $\sc$ and $\sd$
are unital, it is easily seen that these rules yield a well
defined functor $\omega_B:\sc_B\to\sd_B$. To define a coherence
constraint $\tau^\omega$ for this system $(\omega_B~|~B\in\Ob(\cB))$,
notice that for every morphism $\phi:B\to B'$ of $\cB$ and
every $Y\in\Ob(\sc_{B'})$ we have the $\cB$-cartesian morphism
$(\phi,\one_{\sc_\phi Y}):(B,\sc_\phi Y)\to(B',Y)$ of $\cFib(\sc)$,
and let $t:\sd_\phi\circ\omega_{B'}Y\to\omega_B\circ\sc_\phi Y$
be the unique morphism of $\sd_B$ such that
$F(\phi,\one_{\sc_\phi Y})=(\phi,t):
(B,\omega_B\sc_\phi Y)\to(B',\omega_{B'}Y)$; in light of lemma
\ref{lem_distinguished-cleavage}(i) we see that $t$ is an
isomorphism of $\sd_B$, and we set $\tau^\omega_{\phi,Y}:=t^{-1}:
\sd_\phi\circ\omega_{B'}Y\to\omega_B\circ\sc_\phi Y$. To check
the naturality of the rule $Y\mapsto\tau^\omega_{\phi,Y}$, notice
that every morphism $f:X\to Y$ in $\sc_{B'}$ yields a commutative
diagram :
$$
\xymatrix{ (B,\sc_\phi X) \ar[rr]^-{(\one_B,\sc_\phi f)}
\ar[d]_{(\phi,\one_{\sc_\phi X})} & &
(B,\sc_\phi Y) \ar[d]^{(\phi,\one_{\sc_\phi Y})} \\
(B',X) \ar[rr]^-{(\one_{B'},f)} & & (B',Y)
}$$
in $\cFib(\sc)$, whence the identity :
$$
(\phi,\tau^{\omega\ -1}_{\phi,Y})\circ(\one_B,\omega_B\sc_\phi f)=
(\one_{B'},\omega_{B'}f)\circ(\phi,\tau^{\omega\ -1}_{\phi,X})
\qquad
\text{in $\cFib(\sd)$}
$$
which in turn translates as the identity
$\tau^{\omega\ -1}_{\phi,Y}\circ\omega_B\sc_\phi(f)=
\sd_\phi\omega_{B'}(f)\circ\tau^{\omega\ -1}_{\phi,X}$ in $\sd_B$,
as required. Lastly, we need to verify the coherence axioms
for $\tau^\omega$. However, since $\sc$ and $\sd$ are unital,
the first coherence axiom amounts to the identity
$$
\tau^\omega_{\one_B}=\one_{\omega_B}
\qquad
\text{for every $B\in\Ob(\cB)$}
$$
which is immediate from the construction of $\tau^\omega$.
Next, let $\phi:B\to B'$ and $\psi:B'\to B''$ be two morphisms
of $\cB$, and $X$ any object of $\sc_{B''}$; we notice the
commutative diagram :
$$
{\diagram (B,\sc_\phi\sc_\psi X) \ar[d]_{(\one_B,\gamma^\sc_{(\phi,\psi),X})}
\ar[rr]^-{(\phi,\one_{\sc_\phi\sc_\psi X})} & & (B',\sc_\psi X)
\ar[d]^{(\psi,\one_{\sc_\psi X})} \\
(B,\sc_{\psi\phi}X) \ar[rr]^-{(\psi\phi,\one_{\sc_{\psi\phi}X})} & &
(B'',X).
\enddiagram}
\qquad
\text{ in $\cFib(\sc)$}
$$
which, after applying the functor $F$, yields the identity
in $\cFib(\sd)$ :
$$
(\psi,\tau^{\omega\ -1}_{\psi,X})\circ
(\phi,\tau^{\omega\ -1}_{\phi,\sc_\psi X})=
(\psi\phi,\tau^{\omega\ -1}_{\psi\phi,X})\circ
(\one_B,\omega_B(\gamma^\sc_{(\phi,\psi),X})).
$$
which in turn translates as the second coherence axiom.
A direct inspection now gives :
$$
\cFib(\omega)=F.
$$
Summing up, we have shown that \eqref{eq_sPsNat-to-cFib} is
bijective on objects. A simple inspection shows that
\eqref{eq_sPsNat-to-cFib} is also injective on morphisms.
To conclude, consider therefore two pseudo-natural
transformations $\omega,\omega':\sc\to\sd$ and a natural
transformation $\xi:\cFib(\omega)\Rightarrow\cFib(\omega')$
such that $\pi^\sd*\xi=\one_{\pi^\sc}$, and for every object
$(B,X)$ of $\cFib(\sc)$, let $\Xi_{B,X}:\omega_BX\to\omega'_BX$
be the unique morphism of $\sd_B$ such that
$\xi_{(B,X)}=(\one_B,\Xi_{B,X})$. We need to check that the rule:
$(B,X)\mapsto\Xi_{B,X}$ yields a well defined modification
$\omega\leadsto\omega'$, in which case we shall have
$\cFib(\Xi)=\xi$, by inspection. However, for every
$B\in\Ob(\cB)$, the naturality of the rule :
$X\mapsto\Xi_{B,X}$ follows from a simple inspection.
It remains thus only to verify the compatibility condition
for $\Xi$. To this aim, let $\phi:B\to B'$ be any morphism
of $\cB$, and $Y$ any object of $\sc_{B'}$; the morphism
$(\phi,\one_{\sc_\phi Y}):(B,\sc_\phi Y)\to(B',Y)$ of
$\cFib(\sc)$ yields a commutative diagram
$$
{\diagram (B,\omega_B\sc_\phi Y) \ar[rr]^-{\xi_{(B,\sc_\phi Y)}}
\ar[d]_{(\phi,\tau^{\omega\ -1}_{\phi,Y})} & & 
(B,\omega'_B\sc_\phi Y) \ar[d]^{(\phi,\tau^{\omega'\ -1}_{\phi,Y})} \\
(B',\omega_{B'}Y) \ar[rr]^-{\xi_{(B',Y)}} & &
(B',\omega'_{B'}Y)
\enddiagram}
\qquad
\text{in $\cFib(\sd)$
}$$
which in turn is equivalent to the required identity : we
leave the details to the reader.
\end{proof}

\begin{remark}\label{rem_distinguished-cleavage}
(i)\ \
Let $\phi:\cA\to\cB$ be any fibration with small fibres, and
$\sc$ the pseudo-functor associated with a given cleavage
$\blambda$ for $\phi$. Claim \ref{cl_exhibit-nat-isom} exhibits
a natural isomorphism $\psi^\blambda:\cFib(\sc)\isom\cA$
and the fibration $\pi^\sc:\cFib(\sc)\to\cB$ carries
the distinguished cleavage $\blambda^*$ provided
by lemma \ref{lem_distinguished-cleavage}(ii).
By inspecting the constructions, we find
$$
\psi^\blambda(\blambda^*(X,g))=
\psi^\blambda(g,\one_{\sc_\phi X})=\blambda(X,g)
$$
for every morphism $g:B'\to B$ of $\cB$ and every object $X$ of
$\phi^{-1}B$, whence :
$$
\psi^\blambda\circ\blambda^*=\blambda.
$$
Denote by $\sc^*$ the pseudo-functor associated with
$\blambda^*$; it follows that the restrictions of
$\psi^\blambda$
$$
\psi^\blambda_{|B}:\cFib(\sc)_B\isom\cA_B
\qquad
\text{for every $B\in\Ob(\cB)$}
$$
define a strict pseudo-natural isomorphism
of pseudo-functors
$$
\psi^\blambda_{|\bullet}:\sc^*\isom\sc.
$$

(ii)\ \
Let $\sc:\cB^o\to\bCat$ be any pseudo-functor, and
$\rho:\cC\to\cB$ any functor. Then there is a natural
commutative diagram
$$
\xymatrix{ \cFib(\sc\circ\rho^o) \ar[rr]^-{\cFib(\rho)}
\ar[d]_{\pi^{\sc\circ\rho^o}} & & \cFib(\sc) \ar[d]^{\pi^\sc} \\
\cC \ar[rr]^-\rho & & \cB
}$$
where $\cFib(\rho)$ is the functor given by the rules :
$(C,X)\mapsto(\rho(C),X)$ and $(g,f)\mapsto(\rho(g),f)$
for every object $(C,X)$ and morphism $(g,f)$ of
$\cFib(\sc\circ\rho)$. This diagram induces a natural
identification of fibrations over $\cC$ :
$$
\sT^\rho_\sc:\cFib(\sc\circ\rho^o)\isom\Fib(\rho)^*(\cFib(\sc))
\qquad
(C,X)\mapsto(C,(\rho(C),X))
$$
(notation of remark \ref{rem_added-little-extra}(i)).
Especially, via this identification, $\sc\circ\rho^o$ is
the pseudo-functor associated with a cleavage for
$\Fib(\rho)^*(\cFib(\sc))$. By a simple inspection, it is
easily seen that the rule : $\sc\mapsto\sT^\rho_\sc$ yields a
strict pseudo-natural isomorphism of strict pseudo-functors :
$$
\xymatrix{ \sPsFun(\cB^o,\bCat) \ar[rrr]^-{\cFib_\cB}
\ar[d]_{\sPsFun(\rho^o,\bCat)} & \drtwocell\omit{^\ \sT^\rho\ \ } & & 
\Fib(\cB) \ar[d]^{\Fib(\rho)^*} \\
\sPsFun(\cC^o,\bCat) \ar[rrr]_-{\cFib_\cC} & & & \Fib(\cC).
}$$

(iii)\ \
Let $\cA\xrightarrow{\phi}\cB\xleftarrow{\phi'}\cA'$ be two
fibrations with small fibres, $G:\cA\to\cA'$ a $\cB$-cartesian
functor; let us fix cleavages $\blambda$ for $\cA$ and
$\blambda'$ for $\cA'$, and let $\sc$ and $\sc'$ be the
respective associated pseudo-functors. By theorem
\ref{th_fundamental-fibrations}, there exists a
unique pseudo-natural transformation
$$
\sc^G:\sc\Rightarrow\sc'
$$
that makes commute the diagram
$$
\xymatrix{ \cFib(\sc) \ar[r]^-{\psi^\blambda}
\ar[d]_{\cFib(\sc^G)} & \cA \ar[d]^G \\
\cFib(\sc') \ar[r]^-{\psi^{\blambda'}} & \cA'.
}$$
By inspecting the constructions, we can extract the
following explicit description of $\sc^G$ :
\begin{itemize}
\item
for every $B\in\Ob(\cB)$, the functor $\sc^G_B:\cA_B\to\cA'_B$
is the restriction of $G$
\item
for every morphism $g:B\to B'$ in $\cB$, the coherence
constraint of $\sc^G$ is the isomorphism of functors
$\tau_g:\sc'_g\circ\sc^G_B\Rightarrow\sc^G_{B'}\circ\sc_g$
that assigns to every $A\in\Ob(\cA_{B'})$ the unique
isomorphism of $\cA'_B$
$$
\tau_{g,A}:\sc'_g(GA)\isom G(\sc_gA)
\qquad\text{such that}\qquad
G(\blambda(A,g))\circ\tau_{g,A}=\blambda'(GA,g).
$$
\end{itemize}
\end{remark}

\begin{corollary}\label{cor_fibrations}
Let $\phi:\cA\to\cB$ and $\phi':\cA'\to\cB$ be two fibrations,
$F:\cA\to\cA'$ a cartesian $\cB$-functor, and $i\leq 2$ an
integer. We have :
\begin{enumerate}
\item
The following conditions are equivalent :
\begin{enumerate}
\item
$F$ is {\em fibrewise $i$-faithful\/}, {\em i.e.} the
restriction $\cA_B\to\cA'_B$ of the functor $F$ is
$i$-faithful for every $B\in\Ob(\cB)$ (see remark
{\em\ref{rem_i-faithful}}).
\item
$F$ is $i$-faithful, and in case $i=2$, it is even a
$\cB$-equivalence (see remark {\em\ref{rem_added-little-extra}(v)}).
\end{enumerate}
\item
Let $\cC\to\cB$ be any fibration, and suppose that
the conditions of\/ {\em(i)} hold. Then :
\begin{enumerate}
\item
The functor $\sCart_\cB(\cC,F)$ is $i$-faithful.
\item
If\/ $i=2$, also $\sCart_\cB(F,\cC)$ is $i$-faithful.
\end{enumerate}
\end{enumerate}
\end{corollary}
\begin{proof}(i): The implication (i.b)$\Rightarrow$(i.a) is
trivial. For the converse, we consider first the case
where $i=2$ : by theorem \ref{th_fundamental-fibrations}
we may assume that $\cA=\cFib(\sc)$ and $\cA'=\cFib(\sc')$
for (unital) pseudo-functors $\sc,\sc':\cB^o\to\bCat$, and
$F=\cFib(\omega)$ for a pseudo-natural transformation
$\omega:\sc\Rightarrow\sc'$. Then condition (a) means
that $\omega_B:\sc_B\to\sc'_B$ is an equivalence of
categories for every $B\in\Ob(\cB)$. By theorem
\ref{th_pseudo-nat-equiv}, the latter holds if and only
if $\omega$ is a pseudo-natural equivalence, in which
case (i.b) follows.

Next, if $i=0$, let $A,A'\in\Ob(\cA)$ and $f_1,f_2:A\to A'$
be two morphisms of $\cA$ such that $Ff_1=Ff_2$; set
$B:=\phi A$ and $t:=\phi(f_1)$, and notice that
$\phi(f_2)=\phi'(Ff_2)=\phi'(Ff_1)=t$ as well. Pick any
cartesian morphism $h:A''\to A'$ with $\phi(h)=t$; then
there exist unique morphisms $g_1,g_2:A\to A''$ in $\phi^{-1}B$
with $h\circ g_i=f_i$ for $i=1,2$. Therefore
$Fh\circ Fg_1=Fh\circ Fg_2$, and since $F$ is cartesian
we deduce that $Fg_1=Fg_2$; by assumption, it follows that
$g_1=g_2$, whence $f_1=f_2$, as required.

The case where $i=1$ is similar : let $A,A'$ and $B$ as
in the foregoing, and $f':FA\to FA'$ a morphism of $\cA'$;
set $t:=\phi'(f')$, and pick a cartesian morphism $h:A''\to A'$
with $\phi(h)=t$; then $Fh$ is still cartesian, so there
exists a unique morphism $g':FA\to FA''$ in $\phi'^{-1}B$
with $Fh\circ g'=f'$. By assumption, $g'=Fg$ for some
morphism $g:A\to A''$ in $\phi^{-1}B$, so that $F(h\circ g)=f'$,
as required.

(ii): we consider again first the case where $i=2$ : then,
notice that every $\cB$-functor $G:\cA'\to\cA$ that is
quasi-inverse to $F$ is also $\cB$-cartesian (remark
\ref{rem_added-little-extra}(v)); both assertions (ii.a)
and (ii.b) are immediate consequences.

Next, suppose that $i=0$, and let $G,G':\cC\to\cA$ be two
functors, $\alpha,\beta:G\Rightarrow G'$ two natural
transformations with $F*\alpha=F*\beta$. The latter means
that $F(\alpha_C)=F(\beta_C)$ for every $C\in\Ob(\cC)$;
by assumption, we then have $\alpha_C=\beta_C$ for every
such $C$, as required.

Lastly, suppose that $i=1$, and let
$\alpha':F\circ G\Rightarrow F\circ G'$ be a natural
transformation; by assumption, there exists for every
$C\in\Ob(\cC)$ a morphism $\alpha_C:GC\to G'C$ such that
$F(\alpha_C)=\alpha'_C$. We claim that the rule :
$C\mapsto\alpha_C$ yields a natural transformation
$\alpha:G\Rightarrow G'$; indeed, the assertion amounts
to the identity $G'f\circ\alpha_C=\alpha_{C'}\circ Gf$ for
every morphism $f:C\to C'$ of $\cC$. The latter can be
checked after composition with $F$, where it is clear.
\end{proof}

\sset\subsubsection{}\label{subsec_complete-Fib}
By remark \ref{rem_remains-of-the-day}(v), for every universe
$\sV$ with $\sU\subset\sV$, we get a strict pseudo-functor
$$
\sPsFun((-)^o,\sV\tdu\bCat):{}^o\bCat{}^o\to
\sV\tdu\overline{2\tdu\bCat}.
$$
We wish now to show how, likewise, the rule
$\cB\mapsto\sV\tdu\Fib(\cB)$ yields a well defined pseudo-functor.
Indeed, remark \ref{rem_added-little-extra} already says that
every functor $\rho:\cB'\to\cB$ between $\sU$-small categories
induces a strict pseudo-functor
$\sV\tdu\Fib(\rho)^*:\sV\tdu\Fib(\cB)\to\sV\tdu\Fib(\cB')$,
and we have an obvious strict pseudo-natural isomorphism
$\delta^{\sV\tdu\Fib}_\cB:\one_{\sV\tdu\Fib(\cB)}\isom
\sV\tdu\Fib(\one_\cB)^*$ for every small category $\cB$.
Combining with the strict pseudo-natural isomorphisms of
remark \ref{rem_added-little-extra}(vi), we have a datum
$(\delta^{\sV\tdu\Fib}_\bullet,\gamma^{\sV\tdu\Fib}_{\bullet,\bullet})$
fulfilling the unit and associativity axioms for a pseudo-functor
$$
\sV\tdu\Fib:{}^o\bCat^o\to\sV\tdu\overline{2\tdu\bCat}
$$
for which we must however still assign the rule prescribing its
action on natural transformations. To this aim, we choose for
every $\cB\in\Ob(\bCat)$ and every fibration $F:\cA\to\cB$ with
$\sV$-small fibres, a cleavage $\blambda^F$, and we denote by
$\sc^F:\cB^o\to\sV\tdu\bCat$ the corresponding pseudo-functor.
With the notation of remark \ref{rem_distinguished-cleavage},
we deduce for every functor $\rho:\cB'\to\cB$ a cartesian
isomorphism of $\cB'$-fibrations
$$
\sT^\rho_F:=\sV\tdu\Fib(\rho)^*(\psi^{\blambda^F})\circ\sT^\rho_{\sc^F}:
\cFib(\sc^F\circ\rho^o)\isom\sV\tdu\Fib(\rho)^*(\cA).
$$
Let now $\rho,\rho':\cB'\to\cB$ be two functors, and
$\alpha:\rho\Rightarrow\rho'$ a natural transformation.
There exists a unique $\cB'$-cartesian functor
$\sV\tdu\Fib(\alpha)^*_\cA$ making commute the diagram :
$$
\xymatrix@C+30pt{ \cFib(\sc^F\circ\rho'^o) \ar[r]^-{\cFib(\sc^F*\alpha^o)}
\ar[d]_{\sT^{\rho'}_F} & \cFib(\sc^F\circ\rho^o) \ar[d]^{\sT^\rho_F} \\
\sV\tdu\Fib(\rho')^*\cA \ar[r]^-{\sV\tdu\Fib(\alpha)^*_\cA} &
\sV\tdu\Fib(\rho)^*\cA.
}$$
Next, let $F':\cA'\to\cB$ be another fibration, $G:\cA\to\cA'$
a cartesian functor between $\cB$-fibrations, and define
$\sc^G:\sc^F\Rightarrow\sc^{F'}$ as in remark
\ref{rem_distinguished-cleavage}(iii); according to example
\ref{ex_modifications}(ii), the coherence constraint $\tau^G$
of $\sc^G$ defines an invertible modification :
$$
\tau^G_\alpha:(\sc^{F'}*\alpha^o)\odot(\sc^G*\rho'^o)\leadsto
(\sc^G*\rho^o)\odot(\sc^F*\alpha^o).
$$
We consider the induced diagram of $\cB'$-cartesian functors :
$$
\xymatrix@C+30pt{
\cFib(\sc^{F'}\!\circ\rho'^o) \ar[rrr]^-{\cFib(\sc^{F'}*\alpha^o)}
\ar[ddd]_{\sT^{\rho'}_{F'}} & & & \cFib(\sc^{F'}\circ\rho^o)
\ar[ddd]^{\sT^\rho_{F'}} \\
& \ar[lu]^{\cFib(\sc^G*{\rho'^o})} \cFib(\sc^F\!\circ\rho'^o)
\ar[r]^-{\cFib(\sc^F*\alpha^o)} \ar[d]_{\sT^{\rho'}_F} &
\cFib(\sc^F\circ\rho^o) \ar[d]^{\sT^\rho_F}
\ar[ru]_{\cFib(\sc^G*\rho^o)} \\
& \sV\tdu\Fib(\rho')^*\cA \ar[r]^-{\sV\tdu\Fib(\alpha)^*_\cA}
\ar[ld]_{\sV\tdu\Fib(\rho')^*(G)} &
\sV\tdu\Fib(\rho)^*\cA \ar[rd]^{\sV\tdu\Fib(\rho)^*(G)} \\
\sV\tdu\Fib(\rho')^*\cA' \ar[rrr]^-{\sV\tdu\Fib(\alpha)^*_{\cA'}}
& & & \sV\tdu\Fib(\rho)^*\cA'
}$$
whose inner and outer rectangular subdiagrams commute by
definition, and whose left and right trapezoidal subdiagrams
commute by the strict pseudo-naturality of $\sT^\rho$ and
$\sT^{\rho'}$. If we orient the upper trapezoidal subdiagram
with the natural transformation $\cFib(\tau^G_\alpha)^{-1}$,
there exists therefore a unique natural transformation
$$
\tau^\alpha_G:\sV\tdu\Fib(\alpha)^*_{\cA'}\circ\sV\tdu\Fib(\rho')^*(G)
\Rightarrow\sV\tdu\Fib(\rho)^*(G)\circ\sV\tdu\Fib(\alpha)^*_\cA
$$
orienting the lower trapezoidal subdiagram, such that the
resulting oriented cubical diagram (whose other four faces
are oriented by identities) commutes on $2$-cells, in the
sense of \eqref{subsec_transfer-base-change}; the latter
means that :
$$
\sT^\rho_{F'}*\cFib(\tau^G_\alpha)^{-1}=\tau^\alpha_G*\sT^{\rho'}_F.
$$

\begin{lemma} With the notation of \eqref{subsec_complete-Fib},
the rule : $\cA\mapsto\sV\tdu\Fib(\alpha)^*_\cA$ defines a
pseudo-natural transformation $\sV\tdu\Fib(\alpha)^*:
\sV\tdu\Fib(\rho')^*\Rightarrow\sV\tdu\Fib(\rho)^*$
whose coherence constraint is given by the system of natural
transformations $\tau^\alpha_\bullet$.
\end{lemma}
\begin{proof} By inspecting the construction, we are easily
reduced to checking that
$\tau^{G'}_\alpha\boxminus\tau^G_\alpha=\tau^{G'\circ G}_\alpha$
for every pair of cartesian functors of $\cB$-fibrations
$\cA\xrightarrow{G}\cA'\xrightarrow{G'}\cA''$. But this is
clear, since $\sc^{G'\circ G}=\sc^{G'}\odot\sc^G$.
\end{proof}

\sset\subsubsection{}\label{subsec_modify-Fib-alphas}
Next, let $\rho'':\cB'\to\cB$ be another functor, and
$\alpha':\rho'\Rightarrow\rho''$ a second natural transformation.
According to example \ref{ex_modifications}(i), we have an
invertible modification
$$
\gamma^F_{\alpha,\alpha'}:(\sc^F*\alpha^o)\odot(\sc^F*\alpha'^o)
\leadsto\sc^F*(\alpha'\odot\alpha)^o
$$
induced by the coherence constraint $\gamma^F_{\bullet\bullet}$
of $\sc^F$. Whence, a unique natural transformation
$$
\Xi^{\alpha,\alpha'}_\cA:
\sV\tdu\Fib(\alpha)^*_\cA\circ\sV\tdu\Fib(\alpha')^*_\cA
\Rightarrow\sV\tdu\Fib(\alpha'\odot\alpha)^*_\cA
\quad\text{such that}\quad
\sT^\rho_F*\cFib(\gamma^F_{\alpha,\alpha'})=
\Xi^{\alpha,\alpha'}_\cA*\sT^{\rho''}_F.
$$

\begin{lemma}\label{lem_menisco}
With the notation of \eqref{subsec_modify-Fib-alphas}, the following
holds :
\begin{enumerate}
\item
The rule : $\cA\mapsto\Xi^{\alpha,\alpha'}_\cA$ is a
modification $\Xi^{\alpha,\alpha'}\!:\sV\tdu\Fib(\alpha)^*\!\odot\!
\sV\tdu\Fib(\alpha')^*\!\leadsto\!\sV\tdu\Fib(\alpha'\!\odot\!\alpha)^*$.
\item
Especially, we have $\sV\tdu\Fib(\alpha)^*\odot
\sV\tdu\Fib(\alpha')^*=\sV\tdu\Fib(\alpha'\odot\alpha)^*$ in
$\sV\tdu\overline{2\tdu\bCat}$.
\end{enumerate}
\end{lemma}
\begin{proof} Recall that the coherence constraint of
$\sV\tdu\Fib(\alpha)^*\odot\sV\tdu\Fib(\alpha')^*$ assigns to
every cartesian functor $G:\cA\to\cA'$ of $\cB$-fibrations
the natural transformation
$$
\beta:=(\sV\tdu\Fib(\alpha)^*_{\cA'}*\tau^{\alpha'}_G)
\odot(\tau^\alpha_G*\sV\tdu\Fib(\alpha')^*_\cA)
$$
and we need to check the identity :
$$
X:=(\Xi^{\alpha,\alpha'}_{\cA'}*\sV\tdu\Fib(\rho'')^*(G))\odot\beta=
Y:=\tau^{\alpha'\odot\alpha}_G\odot
(\sV\tdu\Fib(\rho)^*(G)*\Xi^{\alpha,\alpha'}_\cA).
$$
To this aim, it suffices to check that $X*\sT^{\rho''}_F=Y*\sT^{\rho''}_F$.
We compute :
$$
\begin{aligned}
X\!\!*\!\!\sT^{\rho''}_F
&\!=\!(\Xi^{\alpha,\alpha'}_{\cA'}*\sV\tdu\Fib(\rho'')^*(G)*\sT^{\rho''}_F)
\odot(\beta*\sT^{\rho''}_F) \\
&\!=\!(\Xi^{\alpha,\alpha'}_{\cA'}\!\!*\!\!\sT^{\rho''}_{F'}\!*\!
\cFib(\sc^G\!*\!\rho''^o)\!)\!\odot\!(\sV\tdu\Fib(\alpha)^*_{\cA'}
\!*\!\!\sT^{\rho'}_{F'}\!*\!\cFib(\tau^G_{\alpha'})^{-1})\!\odot\!
(\tau^\alpha_G\!*\!\sT^{\rho'}_F\!*\!\cFib(\sc^F\!\!*\!\alpha'^o)\!) \\
&\!=\!\sT^\rho_{F'}*(
\cFib(\gamma^{F'}_{\alpha,\alpha'}*(\sc^G*\rho''^o))\odot
\cFib((\sc^{F'}*\alpha^o)*\tau^G_{\alpha'})^{-1}\odot
\cFib(\tau^G_\alpha*(\sc^F*\alpha'^o))^{-1})
\end{aligned}
$$
and :
$$
\begin{aligned}
Y*\sT^{\rho''}_F
&\,=(\sT^\rho_{F'}*\cFib(\tau^G_{\alpha'\odot\alpha})^{-1})\odot
(\sV\tdu\Fib(\rho)^*(G)*\sT^\rho_F*\cFib(\gamma^F_{\alpha,\alpha'})) \\
&\,=(\sT^\rho_{F'}*\cFib(\tau^G_{\alpha'\odot\alpha})^{-1})\odot
(\sT^\rho_{F'}*\cFib((\sc^G*\rho^o)*\gamma^F_{\alpha,\alpha'}).
\end{aligned}
$$
Thus, we are reduced to showing the identity :
$$
\tau^G_{\alpha'\odot\alpha}\odot
(\gamma^{F'}_{\alpha,\alpha'}*(\sc^G*\rho''^o))
=((\sc^G*\rho^o)*\gamma^F_{\alpha,\alpha'})\odot
(\tau^G_\alpha\!*\!(\sc^F\!*\alpha'^o))\odot
((\sc^{F'}\!*\alpha^o)\!*\!\tau^G_{\alpha'}).
$$
which follows directly from the coherence axiom for $\tau^G$.
\end{proof}

In light of lemma \ref{lem_menisco}(ii), the construction of
the pseudo-functor $\sV\tdu\Fib$ of \eqref{subsec_complete-Fib}
shall be complete, once we have shown :

\begin{lemma}
Let $\mu,\mu':\cB''\to\cB'$ and $\rho,\rho':\cB'\to\cB$
be four functors, and $\alpha:\rho\Rightarrow\rho'$,
$\beta:\mu\Rightarrow\mu'$ two natural transformations. Then
the following diagram commutes in $\sV\tdu\overline{2\tdu\bCat}$ :
$$
\xymatrix@C+60pt{ \sV\tdu\Fib(\mu')^*\circ\sV\tdu\Fib(\rho')^*
\ar@{=>}[r]^-{\sV\tdu\Fib(\beta)^**\sV\tdu\Fib(\alpha)^*}
\ar@{=>}[d]_{\gamma^{\sV\tdu\Fib}_{\mu',\rho'}} &
\sV\tdu\Fib(\mu)^*\circ\sV\tdu\Fib(\rho)^*
\ar@{=>}[d]^{\gamma^{\sV\tdu\Fib}_{\mu,\rho}} \\
\sV\tdu\Fib(\rho'\circ\mu')^*
\ar@{=>}[r]^-{\sV\tdu\Fib(\alpha*\beta)^*} &
\sV\tdu\Fib(\rho\circ\mu)^*.
}$$
\end{lemma}
\begin{proof} In view of lemma \ref{lem_menisco}(ii), it suffices
to consider separately the cases where $\mu=\mu'$ and $\beta=\one_\mu$,
and where $\rho=\rho'$ and $\alpha=\one_\rho$. In the first case,
we have to check the identities :
$$
\begin{aligned}
\gamma^{\sV\tdu\Fib}_{(\mu,\rho),\cA}\circ
\sV\tdu\Fib(\mu)^*(\sV\tdu\Fib(\alpha)^*_\cA)&\,=
\sV\tdu\Fib(\alpha*\mu)^*_\cA\circ\gamma^{\sV\tdu\Fib}_{(\mu,\rho'),\cA} \\
X:=\tau^{\alpha*\mu}_G*\gamma^{\sV\tdu\Fib}_{(\mu,\rho'),\cA}&\,=
Y:=\gamma^{\sV\tdu\Fib}_{(\mu,\rho),\cA'}*\sV\tdu\Fib(\mu)^*(\tau^\alpha_G)
\end{aligned}
$$
for every fibration $F:\cA\to\cB$ and every cartesian functor
of $\cB$-fibrations $G:\cA\to\cA'$. For the first identity, we
consider the diagram :
$$\xymatrix@C+20pt{
\cFib(\sc^F\!\circ\rho'^o\circ\mu^o)
\ar[rr]^-{\cFib(\sc^F*\alpha^o*\mu^o)}
\ar[d]_{\sT^\mu_{\sc^F\circ\rho'^o}} & &
\cFib(\sc^F\circ\rho^o\circ\mu^o) \ar[d]^{\sT^\mu_{\sc^F\circ\rho^o}} \\
\sV\tdu\Fib(\mu)^*\cFib(\sc^F\!\circ\rho'^o)
\ar[rr]^-{\sV\tdu\Fib(\mu)^*\cFib(\sc^F*\alpha^o)}
\ar[d]_{\sV\tdu\Fib(\mu)^*(\sT^{\rho'}_F)} & &
\sV\tdu\Fib(\mu)^*\cFib(\sc^F\circ\rho^o)
\ar[d]^{\sV\tdu\Fib(\mu)^*(\sT^\rho_F)} \\
\sV\tdu\Fib(\mu)^*\sV\tdu\Fib(\rho')^*\cA
\ar[rr]^-{\sV\tdu\Fib(\mu)^*\sV\tdu\Fib(\alpha)^*_\cA}
\ar[d]_{\gamma^\Fib_{(\mu,\rho'),\cA}} & &
\sV\tdu\Fib(\mu)^*\sV\tdu\Fib(\rho)^*\cA
\ar[d]^{\gamma^\Fib_{(\mu,\rho),\cA}} \\
\sV\tdu\Fib(\rho'\circ\mu)^*\cA
\ar[rr]^-{\sV\tdu\Fib(\alpha*\mu)^*_\cA} & &
\sV\tdu\Fib(\rho\circ\mu)^*\cA
}$$
whose central square subdiagrams commutes by construction,
and whose upper square subdiagram commutes by strict naturality
of $\sT^\mu$. Notice that the composition of the three left
(resp. right) vertical arrows equals $\sT^{\rho'\circ\mu}_F$ (resp.
$\sT^{\rho\circ\mu}_F$); it then follows that also the external
square diagram commutes. We conclude that the bottom square
subdiagram commutes as well, whence the sought identity.
For the second identity, set
$Z:=\sV\tdu\Fib(\mu)^*(\sT^{\rho'}_F)\circ\sT^\mu_{\sc^F\circ\rho'^o}$;
it suffices to check that $X*Z=Y*Z$. We compute :
$$
\begin{aligned}
X*Z&\,=\tau^{\alpha*\mu}_G*\sT^{\rho'\circ\mu}_F=
\sT^{\rho\circ\mu}_{F'}*\cFib(\tau^G_{\alpha*\mu})^{-1} \\
Y*Z&\,=\gamma^{\sV\tdu\Fib}_{(\mu,\rho),\cA'}*
\sV\tdu\Fib(\mu)^*(\tau^\alpha_G*\sT^{\rho'}_F)*\sT^\mu_{\sc^F\circ\rho'^o} \\
&\,=\gamma^{\sV\tdu\Fib}_{(\mu,\rho),\cA'}*\sV\tdu\Fib(\mu)^*
(\sT^\rho_{F'}*\cFib(\tau^G_\alpha)^{-1})*\sT^\mu_{\sc^F\circ\rho'^o} \\
&\,=\sT^{\rho\circ\mu}_{F'}*\sT^{\mu\,-1}_{\sc^{F'}\circ\rho^o}*
\sV\tdu\Fib(\mu)^*\cFib(\tau^G_\alpha)^{-1}*\sT^\mu_{\sc^F\circ\rho'^o} \\
&\,=\sT^{\rho\circ\mu}_{F'}*\cFib(\tau^G_\alpha*\mu^o)
\end{aligned}
$$
where the last equality holds by the strict pseudo-naturality
of $\sT^\mu$. So we are reduced to checking that
$\tau^G_{\alpha*\mu}=\tau^G_\alpha*\mu^o$. The latter follows by
direct inspection.

For the case where $\alpha=\one_\rho$, let $F:\cA\to\cB$ be any
fibration, set $\rho^*\cA:=\sV\tdu\Fib(\rho)^*\cA$, and let
$\rho^*F:\rho^*\cA\to\cB'$ be the induced functor. There exists a
unique pseudo-natural isomorphism
$$
\omega_{F,\rho}:\sc^F\circ\rho^o\isom\sc^{\rho^*F}
$$
that makes commute the diagram :
$$
\xymatrix{
\cFib(\sc^F\circ\rho^o) \ar[rd]_{\sT^\rho_F}
\ar[rr]^-{\cFib(\omega_{F,\rho})}
& & \cFib(\sc^{\rho^*F}) \ar[ld]^{\psi^{\blambda^{\rho^*F}}} \\
& \rho^*\cA.
}$$
It follows easily that for every cartesian functor $G:\cA\to\cA'$
we have :
\set\begin{equation}\label{eq_variance}
\sc^{\sV\tdu\Fib(\rho)^*G}=
\omega_{F',\rho}\odot(\sc^G*\rho)\odot\omega_{F,\rho}^{-1}
\end{equation}
(details left to the reader). We consider the diagram :
$$
\xymatrix@C-1pt{ \cFib(\sc^{\rho^*F}\!\circ\!\mu'^o)
\ar[rrrr]^-{\cFib(\sc^{\rho^*F}*\beta^o)} \ar[rd]^{\sT^{\mu'}_{\rho^*F}}
\ar[ddd]|{\cFib(\omega_{F,\rho}*\mu'^o)^{-1}} & & & &
\cFib(\sc^{\rho^*F}\!\circ\!\mu^o) \ar[ld]_{\sT^\mu_{\rho^*F}}
\ar[ddd]|{\cFib(\omega_{F,\rho}*\mu^o)^{-1}} \\
& \sV\tdu\Fib(\mu')^*(\rho^*\cA)
\ar[rr]^-{\sV\tdu\Fib(\beta)^*_{\rho^*\cA}}
\ar[d]_{\gamma^{\sV\tdu\Fib}_{(\mu',\rho),\cA}} & &
\sV\tdu\Fib(\mu)^*(\rho^*\cA)
\ar[d]^{\gamma^{\sV\tdu\Fib}_{(\mu,\rho),\cA}} \\
& \sV\tdu\Fib(\rho\circ\mu')^*\cA \ar[rr]^{\sV\tdu\Fib(\rho*\beta)^*_\cA}
& & \sV\tdu\Fib(\rho\circ\mu)^*\cA \\
\cFib(\sc^F\!\circ\!\rho^o\!\circ\!\mu'^o)
\ar[rrrr]^-{\cFib(\sc^F*\rho^o*\beta^o)} \ar[ru]^{\sT^{\rho\circ\mu'}_F}
& & & & \cFib(\sc^F\!\circ\!\rho^o\!\circ\!\mu^o)
\ar[lu]_{\sT^{\rho\circ\mu}_F}
}$$
whose upper and lower trapezoidal subdiagrams commute by
construction; it is easily seen that the same holds for the
right and left trapezoidal subdiagrams (details left to the
reader). Moreover, by example \ref{ex_modifications}(ii) there
exists an invertible modification
$$
\tau^{\omega_{F,\rho}^{-1}}_{\beta^o}:
(\omega_{F,\rho}^{-1}*\mu^o)\odot(\sc^{\rho^*F}*\beta^o)\leadsto
((\sc^F\circ\rho^o)*\beta^o)\odot(\omega_{F,\rho}^{-1}*\mu'^o).
$$
Hence, there exists a unique orientation
$$
\Xi_\cA:\gamma^{\sV\tdu\Fib}_{(\mu,\rho),\cA}\circ
\sV\tdu\Fib(\beta)^*_{\rho^*\cA}\Rightarrow
\sV\tdu\Fib(\beta*\rho)^*_\cA\circ\gamma^{\sV\tdu\Fib}_{(\mu',\rho),\cA}
$$
for the inner square subdiagram, such that :
$$
\Xi_\cA*\sT^{\mu'}_{\rho^*F}=
\sT^{\rho\circ\mu}_F*\cFib(\tau^{\omega_{F,\rho}^{-1}}_{\beta^o}).
$$
To conclude the proof, it suffices then to show :

\begin{claim} The rule : $\cA\mapsto\Xi_\cA$ defines an invertible
modification
$$
\Xi:\gamma^{\sV\tdu\Fib}_{\mu,\rho}\odot
(\sV\tdu\Fib(\beta)^**\sV\tdu\Fib(\rho)^*)\leadsto
\sV\tdu\Fib(\beta*\rho)^*\odot\gamma^{\sV\tdu\Fib}_{\mu',\rho}.
$$
\end{claim}
\begin{pfclaim}[] Notice that the coherence constraint of
$\gamma^{\sV\tdu\Fib}_{\mu,\rho}\odot
(\sV\tdu\Fib(\beta)^**\sV\tdu\Fib(\rho)^*)$ assigns to every
cartesian functor $G:\cA\to\cA'$ of $\cB$-fibrations the
natural transformation
$X:=\gamma^{\sV\tdu\Fib}_{(\mu,\rho),\cA'}*\tau^\beta_{\sV\tdu\Fib(\rho)^*(G)}$,
and the coherence constraint of
$\sV\tdu\Fib(\beta*\rho)^*\odot\gamma^{\sV\tdu\Fib}_{\mu',\rho}$
assigns to $G$ the natural transformation
$Y:=\tau^{\rho*\beta}_G*\gamma^{\sV\tdu\Fib}_{(\mu',\rho),\cA}$. Then
the assertion comes down to the identity:
$$
X':=(\Xi_{\cA'}*(\sV\tdu\Fib(\mu')^*\sV\tdu\Fib(\rho)^*G))\odot X=
Y':=Y\odot((\sV\tdu\Fib(\rho\circ\mu)^*G)*\Xi_\cA)
$$
and it suffices to check that
$X'*\sT^{\mu'}_{\rho^*F}=Y'*\sT^{\mu'}_{\rho^*(F)}$. We compute :
$$
\begin{aligned}
Y'\!*\!\sT^{\mu'}_{\rho^*F}&\,=(\tau^{\rho*\beta}_G*\sT^{\rho\circ\mu'}_F*
\cFib(\omega_{F,\rho}*\mu'^o)^{-1})\odot
((\sV\tdu\Fib(\rho\circ\mu)^*G)*\sT^{\rho\circ\mu}_F*
\cFib(\tau^{\omega_{F,\rho}^{-1}}_{\beta^o})) \\
&\,=(\sT^{\rho\circ\mu}_{F'}*
\cFib(\tau^{G\,-1}_{\rho*\beta}*(\omega_{F,\rho}*\mu'^o)^{-1}))\odot
((\sV\tdu\Fib(\rho\circ\mu)^*G)*\sT^{\rho\circ\mu}_F*
\cFib(\tau^{\omega_{F,\rho}^{-1}}_{\beta^o})) \\
&\,=(\sT^{\rho\circ\mu}_{F'}*
\cFib(\tau^{G\,-1}_{\rho*\beta}*(\omega_{F,\rho}*\mu'^o)^{-1}))\!\odot\!
(\sT^{\rho\circ\mu}_{F'}*\cFib(\sc^G*\rho^o*\mu^o)*
\cFib(\tau^{\omega_{F,\rho}^{-1}}_{\beta^o})) \\
&\,=\sT^{\rho\circ\mu}_{F'}*
\cFib((\tau^{G\,-1}_{\rho*\beta}*(\omega_{F,\rho}*\mu'^o)^{-1})\odot
((\sc^G*\rho^o*\mu^o)*\tau^{\omega_{F,\rho}^{-1}}_{\beta^o})) \\
X*\sT^{\mu'}_{\rho^*F}&\,=
\gamma^{\sV\tdu\Fib}_{(\mu,\rho),\cA'}*\sT^\mu_{\rho^*F'}*
\cFib(\tau^{\sV\tdu\Fib(\rho)^*G}_\beta)^{-1} \\
&\,=\sT^{\rho\circ\mu}_{F'}*\cFib((\omega_{F',\rho}*\mu^o)^{-1}*
(\tau^{\sV\tdu\Fib(\rho)^*G}_\beta)^{-1})
\end{aligned}
$$
and :
$$
\begin{aligned}
\Xi_{\cA'}*(\sV\tdu\Fib(\mu')^*\sV\tdu\Fib(\rho)^*G)*\sT^{\mu'}_{\rho^*F}
&\,=\Xi_{\cA'}*\sT^{\mu'}_{\rho^*F'}*\cFib(\sc^{\sV\tdu\Fib(\rho)^*G}*\mu'^o) \\
&\,=\sT^{\rho\circ\mu}_{F'}*\cFib(\tau_{\beta^o}^{\omega^{-1}_{F',\rho}}*
(\sc^{\sV\tdu\Fib(\rho)^*G}*\mu'^o))
\end{aligned}
$$
so that $X'*Z=
\sT^{\rho\circ\mu}_{F'}*\cFib((\tau_{\beta^o}^{\omega^{-1}_{F',\rho}}*
(\sc^{\sV\tdu\Fib(\rho)^*G}*\mu'^o))\odot
((\omega_{F',\rho}*\mu^o)^{-1}*(\tau^{\sV\tdu\Fib(\rho)^*G}_\beta)^{-1}))$.
We are thus reduced to checking that :
$$
((\sc^G*\rho^o*\mu^o)*\tau^{\omega_{F,\rho}^{-1}}_{\beta^o})
\odot((\omega^{-1}_{F',\rho}*\mu^o)*\tau^{\sV\tdu\Fib(\rho)^*G}_\beta)\!=\!
(\tau^G_{\rho*\beta}*(\omega^{-1}_{F,\rho}*\mu'^o))\odot
(\tau_{\beta^o}^{\omega^{-1}_{F',\rho}}*(\sc^{\sV\tdu\Fib(\rho)^*G}*\mu'^o)).
$$
The latter follows easily from \eqref{eq_variance} : details
left to the reader.
\end{pfclaim}
\end{proof}

\begin{corollary}\label{cor_murder-in-nice}
We have a pseudo-natural equivalence
$$
\cFib_\bullet:\sPsFun((-)^o,\sV\tdu\bCat)\Rightarrow\sV\tdu\Fib
\qquad
\cC\mapsto\cFib_\cC
$$
whose coherence constraint attaches to every functor $\rho$
of small categories the strict pseudo-natural isomorphism
$(\sT^\rho)^{-1}$.
\end{corollary}
\begin{proof} The required coherence axioms for the rule
$\rho\mapsto(\sT^\rho)^{-1}$ hold by direct inspection (details
left to the reader). There remains therefore only to check the
naturality of this rule, relative to natural transformations
$\beta:\rho\Rightarrow\rho'$ between functors
$\rho,\rho':\cB'\to\cB$. The assertion amounts to checking
the following identity in the $2$-category
$\sV\tdu\overline{2\tdu\bCat}$ :
\set\begin{equation}\label{eq_genou-better-today}
(\sT^\rho)^{-1}\odot(\sV\tdu\Fib(\beta)^**\cFib_\cB)=
(\cFib_{\cB'}*\sPsFun(\beta^o,\sV\tdu\bCat))\odot(\sT^{\rho'})^{-1}.
\end{equation}
Thus, let $\sd:\cB^o\to\sV\tdu\bCat$ be any pseudo-functor,
and denote by $F_\sd:\cFib(\sd)\to\cB$ the induced fibration.
With the notation of \eqref{subsec_complete-Fib}, there exists
a unique pseudo-natural isomorphism
$$
\omega^\sd:\sc^{F_\sd}\isom\sd
\qquad\text{such that}\qquad
\cFib(\omega^\sd)=\psi^{\lambda^{F_\sd}}.
$$
We consider the diagram :
$$
\xymatrix@C+10pt{
\cFib(\sc^{F_\sd}\circ\rho'^o) \ar[r]^-{\sT^{\rho'}_{F_\sd}}
\ar[d]_{\cFib(\sc^{F_\sd}*\beta^o)}
& \sV\tdu\Fib(\rho')^*\cFib(\sd) \ar[r]^-{(\sT^{\rho'}_\sd)^{-1}}
\ar[d]|{\sV\tdu\Fib(\beta)^*{\cFib(\sd)}} &
\cFib(\sd\circ\rho'^o) \ar[d]^{\cFib(\sd*\beta^o)} \\
\cFib(\sc^{F_\sd}\circ\rho^o) \ar[r]_-{\sT^\rho_{F_\sd}}
& \sV\tdu\Fib(\rho)^*\cFib(\sd) \ar[r]_-{(\sT^\rho_\sd)^{-1}} &
\cFib(\sd\circ\rho'^o)
}$$
whose left square subdiagram commutes by construction.
Invoking the strict pseudo-naturality of $\sT^\rho$, it
is easily seen that the composition of the two bottom
horizontal arrows agrees with $\cFib(\omega^\sd*\rho^o)$;
likewise, the composition of the two top horizontal
arrows yields $\cFib(\omega^\sd*\rho'^o)$. Now, according
to example \ref{ex_modifications}(ii), there exists an
invertible modification :
$$
\tau^{\omega^\sd}_\beta:(\sd*\beta^o)\odot(\omega^\sd*\rho'^o)\leadsto
(\omega^\sd*\rho^o)\odot(\sc^{F_\sd}*\beta^o).
$$
Hence, there exists a unique isomorphism of functors
$$
\theta^\beta_\sd:\cFib(\sd*\beta^o)\circ(\sT^{\rho'}_\sd)^{-1}
\isom(\sT^\rho_\sd)^{-1}\circ\sV\tdu\Fib(\rho)^*\cFib(\sd)
\quad\text{such that}\quad
\theta^\beta_\sd*\sT^{\rho'}_{F_\sd}=\cFib(\tau^{\omega^\sd}_\beta)
$$
and we are reduced to checking that the rule :
$\sd\mapsto\theta^\beta_\sd$ defines an invertible modification
from the right to the left side of \eqref{eq_genou-better-today}.
Hence, let $\alpha:\sd\Rightarrow\sd'$ be a pseudo-natural
transformation, and define
$\sc^{\cFib(\alpha)}:\sc^{F_\sd}\Rightarrow\sc^{F_{\sd'}}$ as in
remark \ref{rem_distinguished-cleavage}(iii); it is easily
seen that
$$
\sc^{\cFib(\alpha)}=\omega^{\sd'\ -1}\odot\alpha\odot\omega^\sd.
$$
Now, the coherence constraint of the left-hand side of
\eqref{eq_genou-better-today} assigns to $\alpha$ the natural
transformation $X:=(\sT^\rho_{\sd'})^{-1}*\tau^\beta_{\cFib(\alpha)}$,
where $\tau^\beta_{\cFib(\alpha)}$ is characterized by the identity :
$$
\tau^\beta_{\cFib(\alpha)}*\sT^{\rho'}_{F_\sd}=
\sT^\rho_{F_{\sd'}}*\cFib(\tau^{\cFib(\alpha)}_\beta)^{-1}
$$
and $\tau^{\cFib(\alpha)}_\beta:(\sc^{F_{\sd'}}*\beta^o)\odot
(\sc^{\cFib(\alpha)}*\rho'^o)\leadsto
(\sc^{\cFib(\alpha)}*\rho^o)\odot(\sc^{F_\sd}*\beta^o)$ is
given as well by example \ref{ex_modifications}(ii).
According to remark \ref{rem_remains-of-the-day}(ii),
the right-hand side of \eqref{eq_genou-better-today}
assigns to $\alpha$ the natural transformation
$Y:=\cFib(\tau^\alpha_\beta)^{-1}*(\sT^{\rho'}_\sd)^{-1}$, where
again $\tau^\alpha_\beta:(\sd'*\beta^o)\odot(\alpha*\rho'^o)
\Rightarrow(\alpha*\rho^o)\odot(\sd*\beta^o)$ is given
by example \ref{ex_modifications}(ii). We need then to
check that :
$$
Y':=(\theta^\beta_{\sd'}*\sV\tdu\Fib(\rho')^*\cFib(\alpha))\odot Y
=X':=X\odot(\cFib(\alpha*\rho^o)*\theta^\beta_\sd)
$$
and it suffices to check that
$Y'*\sT^{\rho'}_{F_\sd}=X'*\sT^{\rho'}_{F_\sd}$. We compute :
$$
\begin{aligned}
Y'*\sT^{\rho'}_{F_\sd}
&\,=(\theta^\beta_{\sd'}*\sT^{\rho'}_{F_\sd'}*
\cFib(\sc^{\cFib(\alpha)}*\rho'^o))\odot
\cFib(\tau^\alpha_\beta*(\omega^\sd*\rho'^o))^{-1} \\
&\,=(\cFib(\tau^{\omega^{\sd'}}_\beta)*
\cFib(\sc^{\cFib(\alpha)}*\rho'^o))\odot
\cFib(\tau^\alpha_\beta*(\omega^\sd*\rho'^o))^{-1} \\
X'*\sT^{\rho'}_{F_\sd}
&\,=((\sT^\rho_{\sd'})^{-1}*\sT^\rho_{F_{\sd'}}*
\cFib(\tau^{\cFib(\alpha)}_\beta)^{-1})\odot
(\cFib(\alpha*\rho^o)*\cFib(\tau^{\omega^\sd}_\beta)) \\
&\,=(\cFib(\omega^{\sd'}*\rho^o)*\cFib(\tau^{\cFib(\alpha)}_\beta)^{-1})
\odot(\cFib(\alpha*\rho^o)*\cFib(\tau^{\omega^\sd}_\beta))
\end{aligned}
$$
hence we are further reduced to showing that :
$$
((\omega^{\sd'}*\rho^o)*\tau^{\cFib(\alpha)}_\beta)\odot
(\tau^{\omega^{\sd'}}_\beta*(\sc^{\cFib(\alpha)}*\rho'^o))=
((\alpha*\rho^o)*\tau^{\omega^\sd}_\beta)\odot
(\tau^\alpha_\beta*(\omega^\sd*\rho'^o))
$$
but this is clear, since both sides agree with
$\tau^{\alpha\odot\omega^\sd}_\beta$, where $\tau^{\alpha\odot\omega^\sd}$
denotes as usual the coherence constraint of
$\alpha\odot\omega^\sd$.
\end{proof}

\begin{definition}\label{def_split-fibration}
Let $\cB$ be any category, and $\sV$ any universe.

(i)\ \
A {\em split fibration with $\sV$-small fibres\/} over
$\cB$ is a pair $(\phi,\blambda)$ where $\phi:\cA\to\cB$
is a fibration with $\sV$-small fibres and $\blambda$ a
cleavage for $\phi$ whose associated pseudo-functor
$\sc:\cB^o\to\sV\tdu\bCat$ is strict, in which case
we also say that $\blambda$ is {\em split}.

(ii)\ \
For $i=1,2$, let $(\phi_i:\cA_i\to\cB,\blambda_i)$ be two
split fibrations, and $F:\cA_1\to\cA_2$ a cartesian functor.
Let also $\sc_i$ be the pseudo-functor associated with
$\blambda_i$, for $i=1,2$. For every $B\in\Ob(\cB)$, let
$F_B:\phi_1^{-1}B\to\phi_2^{-1}B$ denote the restriction of
$F$. We say that $F$ is {\em a split cartesian functor\/}
$(\phi_1,\sc_1)\to(\phi_2,\sc_2)$ if the rule :
$B\mapsto F_B$ defines a strict pseudo-natural
transformation $\sc_1\Rightarrow\sc_2$.
\end{definition}

\begin{remark}\label{rem_conditions-for-split}
(i)\ \
With the notation of definition \ref{def_split-fibration}(i),
it is easily seen that $\blambda$ is split if and only if for
every $B\in\Ob(\cB)$ we have $\sc_{\one_B}=\one_{\sc_B}$, and for
every pair of morphisms
$B''\xrightarrow{\ h\ }B'\xrightarrow{\ g\ }B$ of $\cB$ we
have $\sc_{g\circ h}=\sc_g\sc_h$ and the following diagram commutes :
$$
{\diagram
\sc_{g\circ h}A \ar[rr]^-{\blambda(\sc_gA,h)}
\ar[rd]_{\blambda(A,g\circ h)} & & \sc_gA \ar[ld]^{\blambda(A,g)} \\
& A
\enddiagram}
\qquad
\text{for every $A\in\Ob(\phi^{-1}B)$}.
$$

(ii)\ \
Likewise, with the notation of definition
\ref{def_split-fibration}(ii), and taking into account
remark \ref{rem_distinguished-cleavage}(iii), we see that
the $\cB$-cartesian functor $F:\cA_1\to\cA_2$ is split
if and only if for every morphism $g:B'\to B$ of $\cB$
and every $A\in\Ob(\phi^{-1}_1B)$ we have
$$
\sc'_g(FA)=F(\sc_gA)
\qquad\text{and}\qquad
F(\blambda(A,g))=\blambda'(FA,g).
$$
(details left to the reader).
\end{remark}

\begin{example}\label{ex_split-fibration}
The fibration $\ss:\cB/F\cC\to\cB$ associated with any functor
$F:\cC\to\cB$ (see example \ref{ex_fibred-cats}(ii)) is split.
Indeed, we have a natural cleavage, defined as follows.
An object of $\cB/\ss(\cB/F\cC)$ is the datum of objects
$B\in\Ob(\cB)$ and $(f:B'\to FC')\in\Ob(\cB/F\cC)$, and
a morphism $g:B\to B'$ in $\cB$; to such a datum we assign
the commutative diagram
$$
\xymatrix{ B \ar[r]^-{f\circ g} \ar[d]_g & FC' \ar[d]^{F\one_{C'}} \\
B' \ar[r]^-f & FC'
}$$
which is an object $\blambda(B,f,g)$ of $\sMorph(\cB/F\cC)$.
Since this object is a $\cB$-cartesian morphism of $\cB/F\cC$,
the rule $(B,f,g)\mapsto\blambda(B,f,g)$ extends uniquely to
a cleavage for $\ss$. By a simple inspection, we find that the
corresponding pseudo-functor $\sc$ associates with every morphism
$g:B\to B'$ of $\cB$ the functor $g/F\cC:B'/F\cC\to B/F\cC$
(notation of \eqref{subsec_fibreovercat}). Then it is clear
that $\sc$ is strict.
\end{example}

\sset\subsubsection{Cartesian sections of a fibration}
\label{subsec_we-mention}
Let $\cB$ be a category and $\sV\subset\sV'$ two universes.
The $\sV'$-small split fibrations over $\cB$ with essentially
$\sV$-small fibres are the objects of a $2$-category
$$
(\sV,\sV')\tdu\Split(\cB)
$$
whose $1$-cells are given by the split cartesian functors,
and whose $2$-cells are the natural $\cB$-transformations.
As usual, we also often write $\sV\tdu\Split(\cB)$ or just
$\Split(\cB)$ for this $2$-category. Furthermore, we have
an obvious forgetful strict pseudo-functor :
$$
\sF:\sV\tdu\Split(\cB)\to\sV\tdu\Fib(\cB)
\qquad
(\phi,\blambda)\mapsto\phi.
$$
On the other hand, if $\cB$ is small, we have as well
a natural strict pseudo-functor
$$
\sC:\Fib(\cB)\to\Split(\cB).
$$
Namely, for a fibration $\phi:\cA\to\cB$ with small fibres,
we consider the strict pseudo-functor
$$
\cA(-):\cB^o\to\bCat
\qquad
B\mapsto\cA(B):=\sCart_\cB(\cB/B,\cA)
$$
that assigns to every morphism $f:B\to B'$ of $\cB$ the functor
$\sCart(f_*,\cA)$ for  (recall that $\cB/B$ is fibred over $\cB$,
by example \ref{ex_fibred-cats}(i); also $f_*:\cB/B\to\cB/B'$ is
the $\cB$-cartesian functor as in \eqref{eq_push-for}).
The objects of $\cA(B)$ are called the {\em cartesian sections}
of $\cA$ over $B$, and the morphisms of $\cA(B)$ are also called
{\em morphisms of cartesian sections}. Then we set
$$
\sC(\cA):=\cFib(\cA(-))
$$
and we let $\sC(\phi):\sC(\cA)\to\cB$ be the resulting
fibration; recall that the fibre category $\sC(\phi)^{-1}B$
is naturally isomorphic to $\cA(B)$ for every $B\in\Ob(\cB)$,
and $\sC(\cA)$ carries a distinguished cleavage whose
associated pseudo-functor corresponds -- under this
identification -- to the strict pseudo-functor $\cA(-)$
(remark \ref{rem_distinguished-cleavage}(i)).

Every $\cB$-cartesian functor $F:\cA_1\to\cA_2$ induces a
strict pseudo-natural transformation
$$
\sCart_\cB(\cB/-,F):\cA_1(-)\Rightarrow\cA_2(-)
\qquad
B\mapsto\sCart_\cB(\cB/B,F)
$$
whence a $\cB$-cartesian functor
$$
\sC(F):=\cFib(\sCart_\cB(\cB/-,F)):
\sC(\cA_1)\to\sC(\cA_2)
$$
and every natural $\cB$-transformation $\alpha:F_1\Rightarrow F_2$
of $\cB$-cartesian functors $F_1,F_2:\cA_1\to\cA_2$ induces a
modification
$$
\sCart_\cB(\cB/-,\alpha):\sCart_\cB(\cB/-,F_1)\leadsto
\sCart_\cB(\cB/-,F_2)
\qquad
B\mapsto\sCart_\cB(\one_{\cB/B},\alpha)
$$
whence a natural $\cB$-transformation
$$
\sC(\alpha)=:\cFib(\sCart_\cB(\cB/-,\alpha)):
\sC(F_1)\Rightarrow\sC(F_2).
$$
Taking into account the discussion of
\eqref{subsec_fib-from-pseudo}-\eqref{subsec_from-mod-to-trans}
we easily conclude that the foregoing rules yield a well defined
strict pseudo-functor $\sC$ as sought.

\sset\subsubsection{The evaluation functor of a fibration}
\label{subsec_eval-functor}
Let $\cB$ be a small category, and $\phi:\cA\to\cB$ a
fibration with small fibres; we claim that there exists
a $\cB$-cartesian functor
$$
\sev^\cA:\sC(\cA)\to\cA
$$
that assigns to every object $(B,G:\cB/B\to\cA)$ of $\sC(\cA)$
its evaluation $G\one_B\in\Ob(\cA_B)$, and to any
morphism
$(f,\alpha:G\Rightarrow G'\circ f_*):(B,G)\to(B',G')$,
the composition
$$
G\one_B\xrightarrow{\ \alpha_{\one_B}\ } G'f
\xrightarrow{\ G'(f/B')\ } G'\one_{B'}
$$
where $f/B':f\to\one_{B'}$ is the morphism of $\cB/B'$
determined by $f$. In order to check that these rules define
a functor, let
$(g,\beta:G'\Rightarrow G''\circ g_*):(B',G')\to(B'',G'')$
be any other morphism; then
$(g,\beta)\circ(f,\alpha)=(g\circ f,(\beta*f_*)\odot\alpha)$,
and the assertion comes down to the identity
$$
G''(g/B'')\circ\beta_{\one_{B'}}\circ G'(f/B')\circ\alpha_{\one_B}=
G''(g\circ f/B'')\circ((\beta*f_*)\odot\alpha)_{\one_B}.
$$
But we have
$$
\beta_{\one_{B'}}\circ G'(f/B')=G''(g_*(f/B'))\circ\beta_f
\quad\text{and}\quad
G''(g/B'')\circ G''(g_*(f/B'))=G''(g\circ f/B'')
$$
whence the contention. Moreover, a simple inspection
shows that every cartesian functor $F:\cA_1\to\cA_2$
of fibred $\cB$-categories yields a commutative diagram
of $\cB$-categories :
$$
\xymatrix{
\sC(\cA_1) \ar[rr]^-{\sC(F)} \ar[d]_{\sev^{\cA_1}} & &
\sC(\cA_2) \ar[d]^{\sev^{\cA_2}} \\
\cA_1 \ar[rr]^-F && \cA_2
}$$
and for every natural $\cB$-transformation
$\alpha:F_1\Rightarrow F_2$ of $\cC$-cartesian functors
$F_1,F_2:\cA_1\to\cA_2$ we have the identity :
$$
\sev^{\cA_2}*\sC(\alpha)=\alpha*\sev^{\cA_1}.
$$
Hence the rule $\cA\mapsto\sev^\cA$ defines a strict
pseudo-natural transformation
$$
\sev^\bullet:\sF\circ\sC\Rightarrow\one_{\Fib(\cB)}
$$

\begin{remark}\label{rem_eval-functor}
With the notation of \eqref{subsec_eval-functor}, let
$G,G'\in\cA(B)$ be two cartesian sections over an object
$B\in\Ob(\cB)$, and $\alpha:G\Rightarrow G'$ a natural
$\cB$-transformation. Then it is easily seen that $\alpha$
is uniquely determined by its evaluation
$\sev^\cA(\one_B,\alpha)=\alpha_{\one_B}:G\one_B\to G'\one_B$.
Indeed, every object $(f:B'\to B)\in\Ob(\cB/B)$ yields a
morphism $f/B:f\to\one_B$ in $\cB/B$, whence a commutative
diagram
$$
\xymatrix{ Gf \ar[r]^-{\alpha_f} \ar[d]_{G(f/B)} &
G'f \ar[d]^{G'(f/B)} \\
G\one_B \ar[r]^-{\alpha_{\one_B}} & G'\one_B
}$$
which determines $\alpha_f$ uniquely, since $G'(f/B)$ is
cartesian. By the same token, given any morphism
$g:G\one_B\to G'\one_B$ in $\phi^{-1}B$, let us define
$\alpha^g_f:Gf\to G'f$ as the unique morphism in $\phi^{-1}B'$
such that $G'(f/B)\circ\alpha^g_f=g\circ G(f/B)$. Then it
is easily seen that the rule $f\mapsto\alpha^g_f$ defines
a natural $\cB$-transformation $\alpha^g:G\Rightarrow G'$
such that $\alpha^g_{\one_B}=g$. If $G'':\cB/B\to\cA$ is another
cartesian functor, and $g':G'\one_B\to G''\one_B$ is any
morphism in $\phi^{-1}B$, clearly we have :
$$
\alpha^{g'}\odot\alpha^g=\alpha^{g'\circ g}.
$$
\end{remark}

\begin{theorem}\label{th_split-fibration}
The pseudo-functor $\sC$ is a strong right $2$-adjoint
for $\sF$, and is fully faithful.
\end{theorem}
\begin{proof} According to proposition
\ref{prop_triangular-identities}(i) and theorem
\ref{th_2-adjunction}(i), in order to prove the first assertion
it suffices to exhibit a unit and a counit fulfilling the triangular
identities \eqref{subsec_adj-pair}. Our candidate counit shall be
the pseudo-natural transformation $\sev^\bullet$. Taking into account
remark \ref{rem_added-little-extra}(v) and corollary
\ref{cor_fully-faith-2-adjoint}, in order to see that $\sev^\cA$ is
$\cB$-cartesian for every fibration $\phi:\cA\to\cB$, and to prove
the second assertion of the theorem, it then suffices to notice :

\begin{claim}\label{cl_split-fibrations}
The functor $\sev^\cA$ is a $\cB$-equivalence.
\end{claim}
\begin{pfclaim} It suffices to show that $\sev^\cA$ is
a fibrewise equivalence (corollary \ref{cor_fibrations}(i)).
However, from remark \ref{rem_eval-functor} we see already
that for every $B\in\Ob(\cB)$ the restriction
$$
\sev^\cA_{|B}:\cA(B)\to\cA_B
$$
of $\sev^\cA$ is fully faithful. To show that $\sev^\cA_{|B}$
is essentially surjective, let $\blambda$ be a unital cleavage
for $\phi$, and $\sc$ its associated unital pseudo-functor. For
every $A\in\Ob(\cA_B)$ consider the functor
$$
i_A:\cB/B\to\cB/\phi\cA
\qquad
(g:B'\to B)\mapsto(A,B,g)
$$
that assigns to every morphism $(h/B):g\to g'$ of $\cB/B$ the
morphism $(A,B,(h/A)):(A,B,g)\to(A,B,g')$ of $\cB/\phi\cA$. Let
also $\ss:\sMorph(\cA)\to\cA$ be the source functor (notation
of \eqref{subsec_Morph-cat}); it suffices to notice that
$$
\beta^\blambda_{B,A}:=\ss\circ\blambda\circ i_A:\cB/B\to\cA.
$$
is a cartesian section with $\beta^\blambda_{B,A}(\one_B)=A$.
\end{pfclaim}

Next, let $(\phi:\cA\to\cB,\blambda)$ be any split fibration with
small fibres, and $\sc$ the strict pseudo-functor associated with
$\blambda$; denote also by $\blambda^*$ the distinguished cleavage
of $\sC(\cA)$, whose associated pseudo-functor is $\cA(-)$ (see
\eqref{subsec_we-mention}). To define a unit for our adjunction,
we need to exhibit a natural split cartesian functor :
$$
(\phi,\blambda)\to(\sC(\cA),\blambda^*).
$$
Now, for every $B\in\Ob(\cB)$ any $A\in\Ob(\cA_B)$ define
$\beta^\blambda_{B,A}\in\cA(B)$ as in the proof of claim
\ref{cl_split-fibrations}; by virtue of remark
\ref{rem_eval-functor}, for every morphism $f:A'\to A$ in
$\cA_B$ there exists a unique natural $\cC$-transformation
$\beta^\blambda_{B,f}:\beta^\blambda_{B,A'}\Rightarrow\beta^\blambda_{B,A}$
such that $(\beta^\blambda_{B,f})_{\one_B}=f$, and the rules
$A\mapsto\beta^\blambda_{B,A}$, $f\mapsto\beta^\blambda_{B,f}$
yield a well defined functor
$$
\beta^\blambda_B:\cA_B\to\cA(B).
$$
According to theorem \ref{th_fundamental-fibrations}, it then
suffices to show that the rule : $B\mapsto\beta^\blambda_B$
defines a strict pseudo-natural transformation
$$
\beta^\blambda:\sc\Rightarrow\cA(-).
$$
The assertion amounts to checking the identities
$$
\beta^\blambda_{B,A}\circ g_*=\beta^\blambda_{B',\sc_gA}
\qquad
\text{for every $(B'\xrightarrow{g}B)\in\Ob(\cB/B)$ and every
$A\in\Ob(\cA_B)$}.
$$
However, we have :
$$
\beta^\blambda_{B',\sc_gA}(g')=\sc_{g'}(\sc_gA)=\sc_{g\circ g'}A=
\beta^\blambda_{B,A}(g\circ g')
\qquad
\text{for every object $g':B''\to B'$ of $\cB/B'$}
$$
which shows that $\beta^\blambda_{B,A}\circ g_*$ and
$\beta^\blambda_{B',\sc_gA}$ agree on objects; to see that they
also agree on morphisms, we consider two objects $f:C\to B'$
and $f':C'\to B'$ of $\cB/B'$ and a morphism $h/B':f\to f'$
of $\cB/B'$ (notation as in the proof of claim
\ref{cl_split-fibrations}), and we notice that, due to
remark \ref{rem_conditions-for-split}(i), the three
subdiagrams of the following diagram commute
$$
\xymatrix{ \sc_{gf}A \ar[rrrr]^-{\beta^\blambda_{B',\sc_gA}(h/B')}
\ar[rrd]_{\blambda(\sc_gA,f)} \ar@/_2.5pc/[rrdd]_{\blambda(A,gf)} & & & &
\ar@/^2.5pc/[lldd]^{\blambda(A,gf')}
\sc_{gf'}A \ar[lld]^{\blambda(\sc_gA,f')} \\
& & \sc_gA \ar[d]^{\blambda(A,g)} \\
& & A }$$
The assertion is an immediate consequence. We need to verify
that the rule $(\phi,\blambda)\mapsto\beta^\blambda$ yields a
strict pseudo-natural transformation
$$
\beta^\bullet:\one_{\Split(\cB)}\Rightarrow\sC\circ\sF.
$$
To this aim, let
$G:(\phi:\cA\to\cB,\blambda)\to(\phi':\cA'\to\cB,\blambda')$
be any split $\cB$-cartesian functor, and let $\sc$ and $\sc'$
be the pseudo-functors associated with $\blambda$ and respectively
$\blambda'$; we come down to showing the commutativity of the
diagram
$$
{\diagram \cA_B \ar[r]^-{\beta^\blambda_B} \ar[d]_{G_B} & \cA(B)
\ar[d]^{\sCart_\cB(\cB/B,G)} \\
\cA'_B \ar[r]^-{\beta^{\blambda'}_B} & \cA'(B)
\enddiagram}
\qquad
\text{for every $B\in\Ob(\cB)$}.
$$
However, the composition
$\sCart_\cB(\cB/B,G)\circ\beta^\blambda_B$ assigns to every
$A\in\Ob(\cA_B)$ the cartesian section $\cB/B\to\cA'$
given by the rule : $(B'\xrightarrow{g}B)\mapsto G(\sc_gA)$,
whereas $\beta^{\blambda'}_B\circ G$ assigns to the same object
the cartesian section given by the rule :
$(B'\xrightarrow{g}B)\mapsto\sc'_g(GA)$; moreover, if
$B''\xrightarrow{g'}B$ is any other object of $\cB/B$ and
$h/B:g\to g'$ is any morphism, we get two commutative
diagrams :
$$
\xymatrix{
G(\sc_gA) \ar[rr]^-{G(\beta^\blambda_{B,A}(h/B))}
\ar[rd]_{G(\blambda(A,g))} & &
G(\sc_{g'}A) \ar[ld]^{G(\blambda(A,g'))} &
\sc'_g(GA) \ar[rr]^-{\beta^{\blambda'}_{B,GA}(h/B)}
\ar[rd]_{\blambda'(GA,g)} & &
\sc'_{g'}(GA) \ar[ld]^{\blambda'(GA,g')} \\
& GA & & & GA.
}$$
Since $G$ is split, remark \ref{rem_conditions-for-split}(ii)
says that these two diagrams coincide, {i.e.} the functors
$\sCart_\cB(\cB/B,G)\circ\beta^\blambda_B$ and $\beta^{\blambda'}_B\circ G$
agree on all objects of $\phi^{-1}B$. To see that they agree
also on morphisms, let $f:A\to A'$ be any morphism of $\phi^{-1}B$;
we have to show that $G*\beta^{\blambda}_{B,f}=\beta^{\blambda'}_{B,Gf}$.
By remark \ref{rem_eval-functor}, it suffices then to compare
the evaluations at $\one_B$ of both of these natural
$\cB$-transformation; but by definition we have
$(G*\beta^{\blambda}_{B,f})_{\one_B}=Gf=(\beta^{\blambda'}_{B,Gf})_{\one_B}$,
as required.

Lastly, let
$G,G':(\phi:\cA\to\cB,\blambda)\to(\phi':\cA'\to\cB,\blambda')$
be two split $\cB$-cartesian functors, and $\gamma:G\Rightarrow G'$
a natural $\cB$-transformation; we need to show :
$$
\sCart_\cB(\cB/B,\gamma)*\beta_B^\blambda=
\beta_B^{\blambda'}*(\gamma_{|\cA_B})
\qquad
\text{for every $B\in\Ob(\cB)$}
$$
which translates as the identity : $\gamma_{\sc_gA}=\sc'_g(\gamma_A)$
for every $g\in\Ob(\cB/B)$ and $A\in\Ob(\cA_B)$. However,
again by virtue of remark \ref{rem_conditions-for-split}(ii), the
following commutative diagrams coincide :
$$
\xymatrix{
G(\sc_gA) \ar[rr]^-{\gamma_{\sc_gA}} \ar[d]_{G(\blambda(A,g))} & &
G'(\sc_gA) \ar[d]^{G'(\blambda(A,g))} & &
\sc'_g(GA) \ar[rr]^-{\sc'_g(\gamma_A)} \ar[d]_{\blambda'(GA,g)} & &
\sc'_g(G'A) \ar[d]^{\blambda'(G'A,g)} \\
GA \ar[rr]^-{\gamma_A} & & G'A & & GA \ar[rr]^-{\gamma_A} & & G'A
}$$
whence the contention.

It remains to show that the unit and counit thus defined
fulfill the triangular identities. Now, let $\phi:\cA\to\cB$
be any fibration, $(\phi':\cA'\to\cB,\blambda')$ any split
fibration, and $\sc'$ the pseudo-functor associated with
$\blambda'$. Denote by $\blambda^*$ the distinguished
cleavage of $\sC(\cA)$ as in \eqref{subsec_we-mention},
whose associated pseudo-functor is (naturally identified
with) $\cA(-)$; then, in view of theorem
\ref{th_fundamental-fibrations}, the triangular identities
come down to the commutativity of the diagrams
$$
\xymatrix{
\cA(B) \ar[r]^-{\beta^{\blambda^*}_B} \ar@{=}[rd] &
\sC(\cA)(B) \ar[d]^{\sCart_\cB(\cB/B,\sev^\cA)} & &
\sc'_B \ar[r]^-{\beta^{\blambda'}_B} \ar@{=}[rd] &
\cA'(B) \ar[d]^{\sev^{\cA'}_{|B}} \\
& \cA(B) & & & \sc'_B
}$$
for every $B\in\Ob(\cB)$. The latter follows by a direct
inspection.
\end{proof}

\begin{remark}\label{rem_mathsfev}
(i)\ \
Let $\phi:\cA\to\cB$ be a fibration, and $\sc$ the
pseudo-functor associated with a given cleavage $\blambda$
for $\phi$. The proof of theorem \ref{th_split-fibration}
furnishes a pseudo-natural equivalence
$$
\sev^\cA_{|\bullet}:\cA(-)\Rightarrow\sc.
$$
In case $\blambda$ is split, it also furnishes a quasi-inverse
for $\sev^\cA_{|\bullet}$ which is a {\em strict} pseudo-natural
equivalence $\beta^\blambda:\sc\Rightarrow\cA(-)$.
However, we cannot deduce from this that the unit of the
$2$-adjunction of theorem \ref{th_split-fibration} is an
equivalence, since we cannot ensure the existence of a
{\em strict} quasi-inverse for $\beta^\blambda$. (Especially,
$\sev^\cA_{|\bullet}$ is not strict, in general.)

(ii)\ \
On the other hand, denote by $\sFib(\cB)$ the full
$2$-subcategory of $\Fib(\cB)$ whose objects are the
fibrations that admit a split cleavage. Then claim
\ref{cl_split-fibrations} implies that the inclusion
functor
$$
\sFib(\cB)\to\Fib(\cB)
$$
is a $2$-equivalence which admits a pseudo-inverse
given by the rule : $(\phi:\cA\to\cB)\mapsto\sC(\cA)$
for every fibration $\phi$.
\end{remark}

\sset\subsubsection{}\label{subsec_connex-comp-fibration}
Notice that the rule that assigns to every small category
$\cC$ its set of connected components $\pi_0(\cC)$ (see
remark \ref{rem_cofinal}(ii)) yields a well defined functor
$$
\pi_0:\bCat\to\Set
$$
that is left adjoint to the fully faithful imbedding
$\Set\to\bCat$ which assigns to every small set the associated
discrete category (see example \ref{ex_universe}(ii)).
Moreover, $\pi_0$ is also a well defined strict pseudo-functor,
for the natural $2$-category structure on $\bCat$ (and here
we regard $\Set$ as a $2$-category whose only $2$-cells are
identities : see example \ref{ex_2-cats}(i)).

Now, if $\cB$ is any category, and $\sc:\cB^o\to\bCat$
any pseudo-functor, it follows that the composition
$\pi_0\circ\sc:\cB^o\to\Set$ is a strict pseudo-functor,
{\em i.e.} it is a presheaf on $\cB$. Also, every
pseudo-natural transformation $\omega:\sc\Rightarrow\sd$
of pseudo-functors $\sc,\sd:\cB^o\to\bCat$ induces a
morphism of presheaves
$\pi_0*\omega:\pi_0\circ\sc\to\pi_0\circ\sd$, and notice
that if $\Xi:\omega\leadsto\omega'$ is any modification
of such pseudo-natural transformations, then
$\pi_0*\omega=\pi_0*\omega'$ and the induced modification
$\pi_0\circ\Xi:\pi_0*\omega\leadsto\pi_0*\omega$ is just the
identity. In view of theorem \ref{th_fundamental-fibrations},
it then follows that the rule $\sc\mapsto\pi_0\circ\sc$ yields
a well defined strict pseudo-functor
$$
\pi^\cB_0:\Fib(\cB)\to\cB^\wedge.
$$
Explicitly, for every fibration $\phi:\cA\to\cB$, the
presheaf $\pi^\cB_0(\cA)$ is given by the rule :
$X\mapsto\pi_0(\phi^{-1}X)$ for every $X\in\Ob(\cB)$.
It is then easily seen that $\pi_0^\cB$ is a strong left
$2$-adjoint to the strict pseudo-functor $\cFib$ of
\eqref{subsec_conversely}.

\begin{remark}
The functor $\cFib$ of \eqref{subsec_conversely} also
admits a right adjoint. Indeed, notice that $\cFib$
factors uniquely as the composition of a strict
pseudo-functor
$$
\cspFib:\cB^\circ\to\Split(\cB)
$$
and the forgetful strict pseudo-functor $\sF$ (notation of
\eqref{subsec_we-mention}). In light of theorem
\ref{th_split-fibration}, it then suffices to find a right
adjoint for $\cspFib$. Now, if $(\phi:\cA\to\cB,\blambda)$
is any split fibration, and $\sc:\cB^o\to\bCat$ the associated
pseudo-functor, notice that the composition
$$
\Ob\circ\sc:\cB^o\to\bCat\to\Set
\qquad
B\mapsto\Ob(\cA_B)
$$
is a presheaf. It is easily seen that the resulting functor
$$
\Ob^\cB:\Split(\cB)\to\cB^o
\qquad
(\phi:\cA\to\cB,\blambda)\mapsto\Ob\circ\sc
$$
is a right adjoint for $\cspFib$. However, we do not get
a right $2$-adjoint by this rule.
\end{remark}

\subsection{\texorpdfstring{$2$}{2}-Fibrations}
\label{sec_2-fibration}
We wish now to consider the generalization of
\eqref{subsec_fib-from-pseudo} where $\cB$ is an arbitrary
$2$-category. In this case, for every pseudo-functor
$\sc:\cB^o\to\bCat$ we get a natural $2$-category structure
on $\cFib(\sc)$, such that the resulting functor
$\pi^c:\cFib(\sc)\to\cB$ is a strict pseudo-functor.
Indeed, for every two objects $(B,X),(B',Y)$ of $\cFib(\sc)$,
and any two morphisms $(\phi_1,f_1),(\phi_2,f_2):(B,X)\to(B',Y)$
of $\cFib(\sc)$, we declare that the $2$-cells
$(\phi_1,f_1)\Rightarrow(\phi_2,f_2)$ are the $2$-cells
$\alpha:\phi_1\Rightarrow\phi_2$ that make commute the
diagram
$$
\xymatrix{ & X \ar[ld]_{f_1} \ar[rd]^{f_2} \\
\sc_{\phi_1}Y \ar[rr]^-{\sc_{\alpha,Y}} & & \sc_{\phi_2}Y.
}$$
If $(\phi_3,f_3):(B,X)\to(B',Y)$ is another morphism
and $\beta:(\phi_2,f_2)\Rightarrow(\phi_3,f_3)$ another
$2$-cell, the composition
$\beta\odot\alpha:(\phi_1,f_1)\Rightarrow(\phi_3,f_3)$
is given by the composition law for $2$-cells of
$\cB(B,B')$. Since
$\sc_{\beta\odot\alpha,Y}=\sc_{\beta,Y}\circ\sc_{\alpha,Y}$, it
is easily seen that this composition rule is well defined.
Likewise, if $(B'',Z)$ is any other object of $\cFib(\sc)$,
$(\psi_1,g_1),(\psi_2,g_2):(B',Y)\to(B'',Z)$ any two morphisms
and $\alpha':(\psi_1,g_1)\Rightarrow(\psi_2,g_2)$ any $2$-cell,
we let $\alpha'*\alpha:(\psi_1,g_1)\circ(\phi_1,f_1)\Rightarrow
(\psi_2,g_2)\circ(\phi_2,f_2)$ be defined by the corresponding
composition law for $2$-cells in $\cB$. In order to check that
this rule is well defined, it suffices to show the commutativity
of the diagram :
$$
\xymatrix{ X \ar[rr]^-{f_1} \ar[rrd]_{f_2} & &
\sc_{\phi_1}Y \ar[rr]^-{\sc_{\phi_1}(g_1)} \ar[d]^{\sc_{\alpha,Y}} & &
\sc_{\phi_1}\sc_{\psi_1}Z \ar[rr]^-{\gamma^\sc_{(\phi_1,\psi_1),Z}}
\ar[d]^{\sc_{\alpha,\sc_{\psi_1}Z}} & &
\sc_{\psi_1\circ\phi_1}Z \ar[dd]^{\sc_{\alpha'*\alpha,Z}} \\
& & \sc_{\phi_2}Y \ar[rr]^-{\sc_{\phi_2}(g_1)}
\ar[rrd]_-{\sc_{\phi_2}(g_2)} & &
\sc_{\phi_2}\sc_{\psi_1}Z \ar[d]^{\sc_{\phi_2}(\sc_{\alpha',Z})} \\
& & & &
\sc_{\phi_2}\sc_{\psi_2}Z \ar[rr]_-{\gamma^\sc_{(\phi_2,\psi_2),Z}}
& & \sc_{\psi_2\circ\phi_2}Z
}$$
(where $\gamma^\sc$ is the coherence constraint of $\sc$); the
latter follows directly from remark \ref{rem_pseudo-funct}(ii)
and the identity :
$\sc_{\phi_2}(\sc_{\alpha',Z})\circ\sc_{\alpha,\sc_{\psi_1}Z}=
(\sc_{\alpha'}*\sc_\alpha)_Z$. It is then obvious that $\cFib(\sc)$
is a $2$-category with these composition laws (and with the
composition law for $1$-cells given by
\eqref{subsec_fib-from-pseudo}), since the required associativity
and unit axioms hold already in $\cB$.

In analogy with the case of usual categories, it is natural
to make the following :

\begin{definition}\label{def_2-fibrations}
Let $\cA$, $\cB$ be two $2$-categories, $\pi:\cA\to\cB$
a strict pseudo-functor.

(i)\ \ 
We say that a $1$-cell $f:A'\to A$ of $\cA$ is
{\em $\cB$-cartesian} if for every $X\in\Ob(\cA)$ the following
diagram is cartesian in the category $\bCat$ :
$$
\xymatrix{
\cA(X,A') \ar[rr]^-{f_*} \ar[d]_\pi & &
\cA(X,A) \ar[d]^\pi \\
\cB(\pi X,\pi A') \ar[rr]^-{(\pi f)_*} & &
\cB(\pi X,\pi A).
}$$

(ii)\ \
We say that $\pi$ (or $\cA$) is a $2$-fibration over $\cB$,
if for every $1$-cell $\phi:B'\to B$ in $\cB$ and every
$A\in\pi^{-1}B$ there exists a $\cB$-cartesian $1$-cell
$f:A'\to A$ in $\cA$ such that $\pi(f)=\phi$.

(iii)\ \
For every $B\in\Ob(\cB)$, the {\em fibre} $\pi^{-1}B$ is
the $2$-subcategory of $\cA$ whose underlying category
is the fibre over $B$ of the functor underlying $\pi$,
and whose $2$-cells are the $2$-cells $\alpha$ of $\cA$
such that $\pi(\alpha)=i_B$. For a universe $\sV$, say
that $\pi$ {\em has small $\sV$-fibres} if $\pi^{-1}B$
is a small $2$-category for every $B\in\Ob(\cB)$.
\end{definition}

\begin{remark} 
(i)\ \
With the notation of definition \ref{def_2-fibrations}, notice
that a $\cB$-cartesian $1$-cell of $\cA$ is also a cartesian
morphism of the underlying category, since the forgetful
functor $\bCat\to\Set$ commutes with fibre products.

(ii)\ \
It is easily seen that a composition of cartesian $1$-cells
is cartesian, and every invertible $1$-cell is cartesian
(details left to the reader).

(iii)\ \
It follows from (i) that if $\pi$ is a $2$-fibration, then
it is also a (usual) fibration on the underlying categories.
Moreover, in this case a $1$-cell of $\cA$ is cartesian if
and only if it is a cartesian morphism of the underlying
category : indeed, the necessity has already been remarked
in (i). Conversely, let $f:A'\to A$ be a $1$-cell which is
a cartesian morphism in the underlying category, and pick
any $1$-cell $g:A''\to A$ of $\cA$ which is cartesian in
the sense of definition \ref{def_2-fibrations}(i), and with
$\pi(g)=\pi(f)$; then there exists an isomorphism
$h:A'\isom A''$ in the underlying category of $\cA$ such that
$f=g\circ h$. By (ii), we deduce that $f$ is a cartesian $1$-cell.
\end{remark}

\begin{lemma}\label{lem_Fib-is-a-2-fibr}
Let $\cB$ be a $2$-category, and $\sc:\cB^o\to\bCat$ a
pseudo-functor. We have :
\begin{enumerate}
\item
The $2$-category $\cFib(\sc)$ is a $2$-fibration over $\cB$.
\item
A $1$-cell $(\phi,f):(B',X)\to(B,Y)$ of $\cFib(\sc)$ is
cartesian if and only if $f:X\to\sc_\phi Y$ is an isomorphism
in $\sc_{B'}$.
\end{enumerate}
\end{lemma}
\begin{proof} It suffices to prove (ii). Now, notice that
the condition of (ii) is necessary, due to lemma
\ref{lem_distinguished-cleavage}(i). Conversely, let
$(\phi,f)$ be such a $1$-cell with $f$ an isomorphism;
consider any object $(B'',Z)$ of $\cFib(\sc)$, any two
$1$-cells $\psi_1,\psi_2:B''\to B'$ in $\cB$, a $2$-cell
$\alpha:\psi_1\Rightarrow\psi_2$ in $\cB$, and any $2$-cell
$\beta:(\phi\circ\psi_1,g_1)\Rightarrow(\phi\circ\psi_2,g_2)$
in $\cFib(\sc)$ such that $\beta=\phi_**\alpha$. We know
already by lemma \ref{lem_distinguished-cleavage}(i)
that for $i=1,2$ there exists a unique $1$-cell
$(\psi_i,g'_i):(B'',Z)\to(B',X)$ of $\cFib(\sc)$ such
that $(\phi,f)\circ(\psi_i,g'_i)=(\phi\circ\psi_i,g_i)$,
and it remains to check that $\alpha$ yields a $2$-cell
$(\psi_1,g'_1)\Rightarrow(\psi_2,g'_2)$. To this aim, we
consider the diagram :
$$
\xymatrix{ Z \ar[rr]^-{g'_1} \ar[rrd]_{g'_2} & &
\sc_{\psi_1}X \ar[d]^-{\sc_{\alpha,X}} \ar[rr]^-{\sc_{\psi_1}(f)} & &
\sc_{\psi_1}\sc_\phi Y \ar[d]^{\sc_{\alpha,\sc_\phi Y}}
\ar[rr]^-{\gamma^\sc_{(\psi_1,\phi),Y}} & &
\sc_{\phi\circ\psi_1}Y \ar[d]^{\sc_{\phi*\alpha,Y}} \\
& & \sc_{\psi_2}X \ar[rr]^-{\sc_{\psi_2}(f)} & & \sc_{\psi_2}\sc_\phi Y
\ar[rr]^-{\gamma^\sc_{(\psi_2,\phi),Y}} & & \sc_{\phi\circ\psi_2}Y
}$$
(where $\gamma^\sc$ is the coherence constraint of $\sc$).
By assumption, the subdiagram obtained after removing the
two middle vertical arrows commutes, and we need to check
that the triangular subdiagram on the left commutes as well.
However, the square subdiagram on the right commutes by
remark \ref{rem_pseudo-funct}(ii), and the central square
subdiagram commutes by naturality oc $\sc_\alpha$. Moreover,
since $f$ is an isomorphism, the two bottom horizontal arrows
are both isomorphisms. The contention follows straightforwardly.
\end{proof}

\begin{example}\label{ex_2-fibrations}
Let $\cB$ be any $2$-category.

(i)\ \
Denote by $\bone$ a chosen final object of $\bCat$, {\em i.e.}
a category with only one object and one morphism, and consider
the constant pseudo-functor $\sF_\bone:\cB^o\to\bCat$ with value
$\bone$ (notation of \eqref{sec_pseudo-cones}). It is easily
seen that the induced projection :
$$
\cFib(\sF_\bone)\to\cB
$$
is an isomorphism of categories.

(ii)\ \
If $\cB$ has small $\Hom$-categories, fix $B\in\cB$
and consider the strict pseudo-functor
$$
h_B:\cB^o\to\bCat
\qquad
B'\mapsto\cB(B',B)
$$
that assigns to every $1$-cell $f:B'\to B''$ the functor
$f^*:\cB(B'',B)\to\cB(B',B)$, and to every $2$-cell
$\beta:f\Rightarrow f'$ between $1$-cells $f,f':B'\to B''$
the natural transformation $\beta^*:f^*\Rightarrow f'^*$
(see example \ref{ex_pseudo-functors}(i)). Then we have
a natural isomorphism of $2$-categories
$$
\cFib(h_B)\isom\cB/B
$$
(notation of example \ref{ex_2-arrows}(ii)) which identifies
the natural projection $\cFib(h_B)\to\cB$ with the source
strict pseudo-functor $\ss_B:\cB/B\to\cB$.
\end{example}

\sset\subsubsection{}\label{subsec_2-Cart}
For any pair of $2$-fibrations $\cA$, $\cA'$ over $\cB$
we denote by
$$
2\tdu\sCart_\cB(\cA,\cA')
$$
the category whose objects are the {\em $\cB$-cartesian
strict pseudo-functors}, {\em i.e.} the strict pseudo-functors
$\cA\to\cA'$ that are $\cB$-cartesian functors on the
underlying categories; the morphisms are the strict
pseudo-natural transformations that are natural
$\cB$-transformations on the underlying categories. Notice
that $2\tdu\sCart_\cB(\cA,\cFib(\sc))$ is a full subcategory
of $\sCart_\cB(\cA,\cFib(\sc))$, for every pseudo-functor
$\sc:\cB^o\to\bCat$ and every $2$-fibration $\cA$.

Just as in definition \ref{def_decartes-functors}(ii),
every pair of $\cB$-cartesian strict pseudo-functors
$H:\cC\to\cA$ and $K:\cA'\to\cC'$ induces a functor
$$
2\tdu\sCart_\cB(H,K):
2\tdu\sCart_\cB(\cA,\cA')\to 2\tdu\sCart_\cB(\cC,\cC')
\qquad
G\mapsto K\circ G\circ H
$$
which assigns to any morphism $\alpha:G\Rightarrow G'$ in
$2\tdu\sCart_\cB(\cA,\cA')$, the natural transformation
$K*\alpha*H:K\circ G\circ H\Rightarrow K\circ G'\circ H$.
As usual, in case $H=\one_\cA$ (resp. $K=\one_{\cA'}$) we also
denote this functor by $2\tdu\sCart_\cB(\cA,K)$ (resp. by
$2\tdu\sCart_\cB(H,\cA')$).

Likewise, if $H,H':\cC\to\cA$ and $K,K':\cA'\to\cC'$ are
four $\cB$-cartesian strict pseudo-functors, every pair of
strict pseudo-natural $\cB$-transformations
$\beta:H\Rightarrow H'$ and $\gamma:K\Rightarrow K'$ induces
a natural transformation :
$$
2\tdu\sCart_\cB(\beta,\gamma):2\tdu\sCart_\cB(H,K)\Rightarrow
2\tdu\sCart_\cB(H',K')
\qquad
G\mapsto\gamma*G*\beta
$$
and again, if $\beta=\one_H$ (resp. $\gamma=\one_K$) we
denote this natural transformation by $2\tdu\sCart_\cB(H,\gamma)$
(resp. $2\tdu\sCart_\cB(\alpha,K)$).

\sset\subsubsection{}
In the situation of \eqref{sec_2-fibration}, let
$\sd:\cB^o\to\bCat$ be another pseudo-functor, and
$\omega:\sc\Rightarrow\sd$ any pseudo-natural transformation.
We claim that the functor
$\cFib(\omega):\cFib(\sc)\to\cFib(\sd)$ defined as in
\eqref{subsec_sd} is a $\cB$-cartesian strict pseudo-functor.
The assertion amounts to saying that any $2$-cell
$\alpha:(\phi_1,f_1)\Rightarrow(\phi_2,f_2)$ as in
\eqref{sec_2-fibration} induces a $2$-cells
$$
\alpha:
\cFib(\omega)(\phi_1,f_1)\Rightarrow\cFib(\omega)(\phi_2,f_2).
$$
In turn, the latter boils down to the commutativity
of the diagram :
$$
\xymatrix{
\omega_B\circ\sc_{\phi_1} Y \ar[d]_{(\tau^\omega_{\phi_1,Y})^{-1}} & &
\ar[ll]_-{\omega_Bf_1} \omega_BX \ar[rr]^-{\omega_Bf_2} & &
\omega_B\circ\sc_{\phi_2}Y \ar[d]^{(\tau^\omega_{\phi_2,Y})^{-1}} \\
\sd_{\phi_1}\circ\omega_{B'}Y \ar[rrrr]^-{\sd_{\alpha,\omega_{B'}Y}}
& & & & \sd_{\phi_2}\circ\omega_{B'}Y.
}$$
However, since $\alpha$ is a $2$-cell in $\cFib(\sc)$, we
have the identity :
$$
\omega_Bf_2=\omega_B(\sc_{\alpha,Y})\circ\omega_Bf_1
$$
and on the other hand, by remark \ref{rem_pseudo-natural}(i),
the naturality of $\tau^\omega$ implies the identity
$$
(\tau^\omega_{\phi_2,Y})^{-1}\circ\omega_B(\sc_{\alpha,Y})=
\sd_{\alpha,\omega_{B'}Y}\circ(\tau^\omega_{\phi_1,Y})^{-1}
$$
whence the assertion. Lastly, for any two pseudo-natural
transformations $\omega,\omega':\sc\Rightarrow\sd$ and
every modification $\Xi:\omega\leadsto\omega'$, a simple
inspection shows that the natural transformation
$\cFib(\Xi):\cFib(\omega)\Rightarrow\cFib(\omega')$
of \eqref{subsec_from-mod-to-trans} is also a strict
pseudo-natural transformation. 

\begin{lemma}\label{lem_already-in-th-fund-fibr}
{\em(i)}\ \
For every pair of pseudo-functors $\sc,\sd:\cB^o\to\bCat$,
the induced functor :
$$
\cFib:\sPsNat(\sc,\sd)\to 2\tdu\sCart_{\cB^o}(\cFib(\sc),\cFib(\sd))
$$
is an isomorphism of categories.

{\em(ii)}\ \
{\em ($2$-Yoneda's Lemma)}\ \
For every $X\in\Ob(\cB)$ and every pseudo-functor
$\sc:\cB^o\to\bCat$ we have an equivalence of categories
pseudo-natural in both arguments :
$$
\sPsNat(h_X,\sc)\isom\sc_X.
$$
\end{lemma}
\begin{proof}(i): Notice that a $\cB$-cartesian strict pseudo-functor
$\cFib(\sc)\to\cFib(\sd)$ is completely determined by its underlying
cartesian functor of (usual) $\cB$-fibrations. Then the assertion
is already known from the proof of theorem
\ref{th_fundamental-fibrations}.

(ii) follows from (i), example \ref{ex_2-fibrations}(ii) and
claim \ref{cl_split-fibrations}.
\end{proof}

\sset\subsubsection{}
Clearly the $2$-fibrations with small fibres over a given
$2$-category $\cB$ form a $2$-category
$$
2\tdu\Fib(\cB)
$$
whose $\Hom$-categories are given by the categories
$2\tdu\sCart_\cB(-,-)$ as in \eqref{subsec_2-Cart}.
We have a strict and strongly faithful pseudo-functor
$$
\cFib_\cB:\sPsFun(\cB^o,\bCat)\to 2\tdu\Fib(\cB)
\qquad
F\mapsto\cFib(F)
$$
but unlike for usual fibrations, $\cFib_\cB$ is not a
$2$-equivalence. As an application of these constructions,
we deduce :

\begin{theorem}\label{th_bCat-cpt}
The $2$-category $\bCat$ is strongly $2$-complete and
strongly $2$-cocomplete.
\end{theorem}
\begin{proof} Let $F:\cB\to\bCat$ be any pseudo-functor
from a small $2$-category $\cB$; we set :
$$
\cL_F:=2\tdu\sCart_{\cB^o}(\cB^o,\cFib(F))
$$
and we define a pseudo-cone $\pi:\sF_{\cL_F}\Rightarrow F$
as follows. First, in light of lemma
\ref{lem_already-in-th-fund-fibr}(i) and example
\ref{ex_2-fibrations}, we have a natural isomorphism
of categories
$$
\omega:\cL_F\isom\sPsNat(\sF_\one,F)
$$
where $\sF_\one:\cB\to\bCat$ is the constant pseudo-functor
with value $\one$. Now, let $\beta:\sF_\one\Rightarrow F$ be any
pseudo-cone; for every $B\in\Ob(\cB)$ we have then the functor
$\beta_B:\one\to FB$, which is identified with the object
$\beta_B(\emptyset)\in\Ob(FB)$, under the isomorphism of
categories :
$$
\bFun(\one,FB)\isom FB
\qquad
G\mapsto G(\emptyset)
\qquad
(\alpha:G\Rightarrow G')\mapsto
\alpha_\emptyset: G(\emptyset)\to G'(\emptyset)
$$
(here $\Ob(\one)=\{\emptyset\}$). Clearly, the rule
$\beta\mapsto\beta_B$ yields a functor
$$
\pi_B:\sPsNat(\sF_\one,F)\to FB
$$
that assigns to every modification $\Xi:\beta\leadsto\beta'$
the natural transformation $\Xi_B:\beta_B\Rightarrow\beta'_B$,
identified with the morphism
$\Xi_{B,\emptyset}:\beta_B(\emptyset)\to\beta'_B(\emptyset)$ in
$FB$. Lastly, it is easily seen that the rule
$B\mapsto\pi_B\circ\omega$ yields a pseudo-cone as sought,
with coherence constraint
$$
\tau^\pi_f:Ff\circ\pi_B\Rightarrow\pi_{B'}
\qquad\text{such that}\qquad
\tau^\pi_{f,\beta}:=\tau^\beta_{f,\emptyset}
$$
for every pseudo-cone $\beta:\sF_\one\to F$ with coherence
constraint $\tau^\beta$.

Explicitly, for every $B\in\Ob(\cB)$
the functor $\pi_B:\cL_F\to FB$ is determined by the identity :
$$
\phi(B)=(B,\pi_B(\phi))
\qquad
\text{for every cartesian strict pseudo-functor
$\phi:\cB^o\to\cFib(F)$}
$$
and for every pair of cartesian strict pseudo-functors
$\phi,\phi':\cB^o\to\cFib(F)$ and every strict pseudo-natural
$\cB^o$-transformation $\beta:\phi\Rightarrow\phi'$, we have
$\beta_B=(\one_B,\pi_B(\beta))$.
 
We claim that $(\cL_F,\pi)$ is a strong $2$-limit of $F$,
{\em i.e.} for every small category $\cC$, the functor
\set\begin{equation}\label{eq_iso-not-just-eq}
\bFun(\cC,\cL_F)\to\sPsNat(\sF_\cC,F)
\qquad
(\phi:\cC\to\cL_F)\mapsto\pi\odot\sF_\phi
\end{equation}
is an isomorphism of categories (see definition
\ref{def_pseudo-lim}(i)). Indeed, let
$\alpha:\sF_\cC\Rightarrow F$ be any pseudo-cone with vertex
$\cC$ and basis $F$; there follows a $\cB^o$-cartesian strict
pseudo-functor $\cFib(\alpha):\cFib(\sF_\cC)\to\cFib(F)$,
whence a functor
$$
\alpha^\dagger:=2\tdu\sCart_{\cB^o}(\cB^o,\cFib(\alpha)):
2\tdu\sCart_{\cB^o}(\cB^o,\cFib(\sF_\cC))\to\cL_F.
$$
If $\Theta:\alpha\leadsto\beta$ is any modification,
we set as well
$\Theta^\dagger:=2\tdu\sCart_{\cB^o}(\cB^o,\cFib(\Theta)):
\alpha^\dagger\to\beta^\dagger$. However, a simple inspection
shows that $\cFib(\sF_\cC)=\cB^o\times\cC$, the product in
the category of $2$-categories (where $\cC$ is regarded as
a $2$-category with trivial $2$-cells), fibered over $\cB^o$
via the natural projection $p:\cB^o\times\cC\to\cB^o$. There
is an obvious functor
$$
c^\cC:\cC\to 2\tdu\sCart_{\cB^o}(\cB^o,\cFib(\sF_\cC))
\qquad
X\mapsto c_X
$$
that assigns to every object $X$ of $\cC$ the constant functor
$c_X:\cB^o\to\cC$ with value $X$, naturally identified with a
section $\cB^o\to\cB^o\times\cC$ of the projection $p$. Then,
a simple inspection shows that the functor
$$
\sPsNat(\sF_\cC,F)\to\bFun(\cC,\cL_F)
\qquad
\alpha\mapsto \alpha^\dagger\circ c^\cC
\qquad
\Theta\mapsto \Theta^\dagger*c^\cC
$$
is inverse to the functor \eqref{eq_iso-not-just-eq}, whence
the contention.

Next we consider the category $\bar\cFib(F)$ such that
$\Ob(\bar\cFib(F)):=\Ob(\cFib(F))$ and
$$
\Hom_{\bar\cFib(F)}(X,X'):=\pi_0(\cFib(F)(X,X'))
\qquad
\text{for every $X,X'\in\Ob(\cFib(F))$}
$$
(notation of \eqref{subsec_connex-comp-fibration}). The
composition law for morphisms in $\bar\cFib(F)$ is induced
in the obvious way by that of $1$-cells of $\cFib(F)$.
For every $1$-cell $f$ of $\cFib(F)$ we denote by $[f]$
the class of $f$ in $\rMorph(\bar\cFib(F))$ and we let
$\Sigma:=\{[f]\in\rMorph(\bar\cFib(F))~|~
\text{$f$ is $\cB$-cartesian}\}$. We set
$$
\cL'_F:=\bar\cFib(F))[\Sigma^{-1}]
$$
(notation of theorem \ref{th_localize-cats}). We have a
natural functor
$$
\cFib(F)\to\cL'_F
$$
namely, the composition of the obvious projection
$\cFib(F)\to\bar\cFib(F)$ with the localization
$\bar\cFib(F)\to\cL'_F$. Recall that $\cFib(F)$ carries
a distinguished cleavage $\blambda^*$, whose associated
pseudo-functor is naturally identified with $F$ (remark
\ref{rem_distinguished-cleavage}(i)). Then, the pseudo-cocone
$\iota:F\to\sF_{\cFib(F)}$ associated with $\blambda^*$ induces
a a pseudo-cocone $\iota':F\Rightarrow\sF_{\cL'_F}$. We claim
that $(\cL'_F,\iota')$ is a strong $2$-colimit of $F$, {\em i.e.}
for every small category $\cC$, the induced functor
\set\begin{equation}\label{eq_colim-iso-not-just-eq}
\bFun(\cL'_F,\cC)\to\sPsNat(F,\sF_\cC)
\qquad
(\phi:\cL'_F\to\cC)\mapsto\sF_\phi\odot\iota'
\end{equation}
is an isomorphism of categories. Indeed, let
$\alpha:F\Rightarrow\sF_\cC$ be any pseudo-cocone; there
follows a $\cB^o$-cartesian strict pseudo-functor
$\cFib(\alpha):\cFib(F)\to\cFib(\sF_\cC)\isom\cB^o\times\cC$.
After composing with the projection $p:\cB^o\times\cC\to\cC$
we get a functor
$$
\alpha^*:\cFib(F)\to\cC.
$$
Since $\cC$ has trivial $2$-cells, $\alpha^*$ factors through
a unique functor $\bar\alpha{}^*:\bar\cFib(F)\to\cC$. By lemma
\ref{lem_Fib-is-a-2-fibr}(ii), a $1$-cell $(f,g):(B,X)\to(B',X')$
of $\cFib(\sF_\cC)$ is cartesian if and only if $g:X\to X'$ is
an isomorphism of $\cC$. It follows that $\bar\alpha{}^*$ in
turn factors uniquely through a functor
$$
\alpha^\ddagger:\cL'_F\to\cC.
$$
On the other hand, notice that a functor $\cFib(F)\to\cC$
factors through $\bar\cFib(F)$ if and only if it induces
a strict pseudo-functor $\cFib(F)\to\cB^o\times\cC$, and
the latter is $\cB^o$-cartesian if and only if the resulting
functor $\bar\cFib(F)\to\cC$ factors through $\cL'_F$.

If $\Theta:\alpha\leadsto\beta$ is any modification, we
get a natural $\cB$-transformation
$\Theta^*:=p*\cFib(\Theta):\alpha^*\Rightarrow\beta^*$,
which in turn induces a natural transformation
$$
\Theta^\ddagger:\alpha^\ddagger\Rightarrow\beta^\ddagger.
$$
By a direct inspection, we see that the functor
$$
\sPsNat(F,\sF_\cC)\to\bFun(\cL'_F,\cC)
\qquad
\alpha\mapsto\alpha^\ddagger
\qquad
\Theta\mapsto\Theta^\ddagger
$$
is inverse to \eqref{eq_colim-iso-not-just-eq}, whence the
contention.
\end{proof}

\begin{example}\label{ex_2-products}
(i)\ \
For instance, let $\cC_\bullet:=(\cC_i~|~i\in I)$ be any
small family of small categories. By inspecting the proof of
theorem \ref{th_bCat-cpt}, we see that the $2$-product
of $\cC_\bullet$ is represented by the product
$\cC:=\prod_{i\in I}\cC_i$, as constructed in example
\ref{ex_cat-cats}(i), and the universal cone for $\cC$
is also a universal pseudo-cone for the $2$-product
of $\cC_\bullet$.

(ii)\ \
Let $F:\cC\to\cB$ and $F':\cC'\to\cB$ be two functors between
small categories; by inspecting the proof of theorem
\ref{th_bCat-cpt}, we see that the $2$-fibre product of $F$
and $F'$ (remark \ref{rem_pseudo-limit}(v)) is represented by
the category whose objects are all data of the form
$X:=(c,c',\xi)$, where $c\in\Ob(\cC)$, $c'\in\Ob(\cC')$, and
$\xi:F(c)\isom F'(c')$ is an isomorphism in $\cB$. If
$X':=(d,d',\zeta)$ is another such datum, the morphisms
$X\to X'$ are the pairs $(\phi,\phi')$, where $\phi:c\to d$
(resp. $\phi':c'\to d'$) is a morphism in $\cC$ (resp. in $\cC'$),
and $\zeta\circ F(\phi)=F'(\phi')\circ\xi$.
\end{example}

\begin{example}\label{ex_filter-2-colim-in-Cat}
(i)\ \
Consider a small category $\cB$ and a pseudo-functor
$F:\cB\to\bCat$. In this case, obviously $\bar\cFib(F)=\cFib(F)$,
so the strong $2$-colimit of $F$ is represented by
$\cFib(F)[\Sigma^{-1}]$, where $\Sigma\subset\rMorph(\cFib(F))$ is
the set of cartesian morphisms.

(ii)\ \
Suppose next that $\cB$ is {\em filtered}. Then we claim that
$\Sigma$ admits a right calculus of fractions. Indeed, it is
clear that $\Sigma$ satisfies conditions (CF1) and (CF2) of
definition \ref{def_right-calculus}(i). Next, let
$(\phi,f):(B,X)\to(B',Y)$ and $(\psi,g):(B'',Z)\to(B',Y)$
be two morphisms of $\cFib(F)$, with $(\psi,g)\in\Sigma$;
hence $\phi:B'\to B$ and $\psi:B'\to B''$ are morphisms of
$\cB$, so there exist morphisms $\phi':B\to B'''$ and
$\psi':B''\to B'''$ in $\cB$ with
$\phi'\circ\phi=\psi'\circ\psi$. It follows that
$(\phi',\one_{F_{\phi'}X}):(B''',F_{\phi'}X)\to(B,X)$ lies in
$\Sigma$, and since $(\psi,g)$ is cartesian there exists
a unique morphism $(\psi',h):(B''',F_{\phi'}X)\to(B'',Y)$ such
that $(\psi,g)\circ(\psi',h)=(\phi,f)\circ(\phi',\one_{F_{\phi'}X})$.
This proves that condition (CF3) holds as well for $\Sigma$.
Lastly, let $(\phi,f),(\phi',f'):(B,X)\to(B',Y)$ and
$(\psi,g):(B',Y)\to(B'',Z)$ be three morphisms of $\cFib(F)$,
with $(\psi,g)$ cartesian, and suppose that
$(\psi,g)\circ(\phi,f)=(\psi,g)\circ(\phi',f')$; hence
$\phi,\phi':B'\to B$ are two morphisms of $\cB$, and by
assumption there exists a morphism $\mu:B\to B'''$ such
that $\mu\circ\phi=\mu\circ\phi'$. Then we have the cartesian
morphism $(\mu,\one_{F_\mu X}):(B''',F_\mu X)\to(B,X)$, and
obviously $(\psi,g)\circ(\phi,f)\circ(\mu,\one_{F_\mu X})=
(\psi,g)\circ(\phi',f')\circ(\mu,\one_{F_\mu X})$. Since
$(\psi,g)$ is cartesian, it follows that
$(\phi,f)\circ(\mu,\one_{F_\mu X})=(\phi',f')\circ(\mu,\one_{F_\mu X})$.
This shows that condition (CF4) holds for $\Sigma$.

(iii)\ \
In particular, in the situation of (ii), the morphisms
in $\cFib(F)[\Sigma^{-1}]$ can be expressed as fractions
with denominators in $\Sigma$, as detailed by proposition
\ref{prop_calculus-frac} (and remark \ref{rem_left-calc-fract}).

(iv)\ \
Let $I$ be a small filtered category, and $F:I\to\bCat$ any
functor; let $\cL$ be the colimit of $F$, as explicitly
given by example \ref{ex_fil-colim-in-Cat}. In particular,
the objects of $\cL$ are equivalence classes $[i,X]$ of pairs
$(i,X)\in\Ob(\cFib(F))$, and the morphisms $[f]:[i,X]\to[i',X']$
are equivalence classes of morphisms $f:F_\phi X\to F_{\phi'}X'$,
for all pairs of morphism
$(i\xrightarrow{\phi}j\xleftarrow{\phi_2}i')$ of $I$. By direct
inspection of the construction of $\cL$, we get a well defined
functor
$$
\cFib(F)\to\cL
\qquad
(i,X)\mapsto[i,X]
\qquad
((i'\to i,f):(i,X)\to(i',X'))\mapsto[f].
$$
By the universal property of the localization, this functor
extends uniquely to a well defined functor
$\cFib(F)[\Sigma^{-1}]\to\cL$, where $\Sigma$ is as in (ii), and
in view of (iii) and lemma \ref{lem_distinguished-cleavage}(i),
it is easily seen that the latter functor is an equivalence. In
other words, the colimit of $F$ represents as well the $2$-colimit
of the same functor, and more precisely, every universal cocone
$F\Rightarrow c_\cL$ is also a (strict) universal pseudo-cocone.
\end{example}

\sset\subsubsection{}\label{subsec_adjoints-of-Fib(rho)}
Let $\sV$ be a universe, $\cC$ a $\sV$-small category, $\cB$
a category with $\sV$-small $\Hom$-sets, and $\rho:\cC\to\cB$
a functor; theorem \ref{th_bCat-cpt} enables us to construct
both right and left $2$-adjoints for the strict pseudo-functor
$\sV\tdu\Fib(\rho)^*$ of remark \ref{rem_added-little-extra}(i).
Indeed, combining with theorem \ref{th_2-Kan-extensions} and
remark \ref{rem_strong-2-Kan-ext}(ii) we see first that the
strong left $2$-Kan extension along $\rho^o$ is a strict left
$2$-adjoint $\sL$ for $\sP:=\sPsFun(\rho^o,\sV\tdu\bCat)$. On
the other hand, for every category $\cB$, by remark
\ref{rem_strict-strong}, the strict and strong $2$-equivalence
$\cFib_\cB$ of theorem \ref{th_fundamental-fibrations} admits
a strict and strong pseudo-inverse
$$
\sQ_\cB:\sV\tdu\Fib(\cB)\isom\sPsFun(\cB^o,\sV\tdu\bCat).
$$
It follows that the strict pseudo-functor
$$
\sV\tdu\Fib(\rho)_!:=\cFib_\cB\circ\sL\circ\sQ_\cC:
\sV\tdu\Fib(\cC)\to\sV\tdu\Fib(\cB)
$$
is a left $2$-adjoint for $\cFib_\cC\circ\sP\circ\sQ_\cB$.
On the other hand, remark \ref{rem_distinguished-cleavage}(ii)
yields a pseudo-natural isomorphism
$\cFib_\cC\circ\sP\isom\sV\tdu\Fib(\rho)^*\circ\cFib_\cB$,
so $\sV\tdu\Fib(\rho)_!$ is also left $2$-adjoint to
$\sV\tdu\Fib(\rho)^*$. Likewise, from the strong right
$2$-Kan extension $\sR$ along $\rho^o$ we get a strict
right $2$-adjoint for $\sV\tdu\Fib(\rho)^*$, namely the
pseudo-functor
$$
\sV\tdu\Fib(\rho)_*:=\cFib_\cB\circ\sR\circ\sQ_\cC:
\sV\tdu\Fib(\cC)\to\sV\tdu\Fib(\cB).
$$
As usual, we will omit mentioning $\sV$, when no ambiguities
are likely to arise from the omission.

\begin{remark}\label{rem_explicit-Fib_*}
(i)\ \
Like all $2$-adjoints, the pseudo-functors $\Fib(\rho)_!$
and $\Fib(\rho)_*$ are well defined up to pseudo-natural
equivalence. However, our construction yields more canonical
representatives, that are well defined up to pseudo-natural
{\em isomorphisms}, since it relies on the strong $2$-Kan
extensions along $\rho^o$. Moreover, both these $2$-adjoints
factor naturally through pseudo-functors
$$
\Fib(\cC)\to\Split(\cB)
$$
again because this is a feature of strong $2$-Kan extensions :
see remark \ref{rem_strong-2-Kan-ext}(ii).

(ii)\ \
In the situation of \eqref{subsec_adjoints-of-Fib(rho)},
we can describe more explicitly the pseudo-functor
$\Fib(\rho)_*$ as follows. Let $\cA\to\cC$ be any fibration
with small fibres; then $\sC(\Fib(\rho)_*(\cA))$ is the
fibration associated with the strict pseudo-functor
$\Fib(\rho)_*(\cA)(-):\cB^o\to\bCat$, and by example
\ref{ex_fibred-cats-II}(ii), for every $B\in\Ob(\cB)$
we have an equivalence of categories :
$$
\Fib(\rho)_*(\cA)(B)\isom\sCart_\cC(\rho\cC/B,\cA)
$$
that is pseudo-natural with respect to morphisms $B'\to B$
in $\cB$. Thus, $\Fib(\rho)_*(\cA)$ is naturally equivalent
to the fibration associated with the strict pseudo-functor
$$
\sCart_\cC(\rho\cC/-,\cA):\cB^o\to\bCat
\qquad
B\mapsto\sCart_\cC(\rho\cC/B,\cA)
$$
that assigns to every morphism $f:B'\to B$ of $\cB$ the
functor $\sCart_\cC(\rho\cC/f,\cA)$.

(iii)\ \
Moreover, we have a pseudo-commutative diagram of
$2$-categories :
$$
\xymatrix{ \cC^\wedge \ar[rr]^-{\cFib_\cC} \ar[d]_{\rho_*} & &
\Fib(\cC) \ar[d]^{\Fib(\rho)_*} \\
\cB^\wedge \ar[rr]^-{\cFib_\cB} & & \Fib(\cB)
}$$
(see remark \ref{rem_pseudo-commutative}). Indeed, let
$j:\Set\to\bCat$ be the inclusion functor (that assigns to
every set the associated discrete category); by inspecting
the constructions of $\Fib(\rho)_*$ and $\rho_*$ (see remark
\ref{rem_was-cofinal}(i)), we come down to exhibiting an
equivalence :
$$
\qquad\qquad
2\tdu\!\!\!\int_{\rho^o} j\circ F\isom j\circ\!\!\int_{\rho^o} F
\qquad\qquad
\text{for every presheaf $F:\cC^o\to\Set$}
$$
pseudo-natural with respect to $F$. Now, for every $B\in\Ob(\cB)$
the category $2\tdu\!\!\int_{\rho^o}j\circ F(B)$ is the $2$-limit
of the functor $j\circ F\circ\st_{B^o}:B^o/\phi^o\cC^o\to\bCat$.
But since $j$ is fully faithful and admits a left $2$-adjoint
(see \eqref{subsec_connex-comp-fibration}), proposition
\ref{prop_2-cat-mitchell} says that such $2$-limit is the
discrete category associated with the set representing the
$2$-limit of $F\circ\st_{B^o}:B^o/\phi^o\cC^o\to\Set$. By remark
\ref{rem_pseudo-limit}(vii), the latter is represented by the
(usual) limit of the same functor. But this in turn is none else
than the definition of $\int_\phi F(B^o)$, whence the contention.

(iv)\ \
In view of corollary \ref{cor_2-Kan-ext-of-fully-faith}, we
also see that if $\rho$ is a fully faithful functor, then
$\Fib(\rho)_*$ and $\Fib(\rho)_!$ are both fully faithful
pseudo-functors.

(v)\ \
For every other universe $\sV'$ with $\sV\subset\sV'$,
we get a diagram of pseudo-functors :
$$
\xymatrix@C+20pt{
\sV\tdu\Fib(\cC) \ar[r]^-{\sV\tdu\Fib(\rho)_!} \ar[d] &
\sV\tdu\Fib(\cB) \ar[d] &
\sV\tdu\Fib(\cC) \ar[r]^-{\sV\tdu\Fib(\rho)_*} \ar[d] &
\sV\tdu\Fib(\cB) \ar[d] \\
\sV'\tdu\Fib(\cC) \ar[r]^-{\sV'\tdu\Fib(\rho)_!} &
\sV'\tdu\Fib(\cB) &
\sV'\tdu\Fib(\cC) \ar[r]^-{\sV'\tdu\Fib(\rho)_*} &
\sV'\tdu\Fib(\cB)
}$$
whose vertical arrows are the inclusions.
A direct inspection of the constructions easily shows
that both diagrams are commutative : the details shall
be left to the reader.

(vi)\ \
Let $\cC$ and $\cC'$ be two small categories, $u:\cC\to\cC'$
and $v:\cC'\to\cC$ two functors, such that $v$ is left adjoint
to $u$. In light of \eqref{subsec_comma-isoms-from-adj}, we
get for every fibration $\cA\to\cC'$ and every $X\in\Ob(\cC)$
a natural isomorphism of categories :
$$
\omega_X:\sCart_{\cC'}(v\cC'/X,\cA)\isom\cA(uX)
$$
and for every morphism $f:Y\to X$ in $\cC$, a commutative
diagram of categories :
$$
\xymatrix{ \sCart_{\cC'}(v\cC'/X,\cA) \ar[r]^-{\omega_X}
\ar[d]_{\sCart_{\cC'}(v\cC'/f,\cA)} & \cA(uX) \ar[d]^{\cA(uf)} \\
\sCart_{\cC'}(v\cC'/Y,\cA) \ar[r]^-{\omega_Y} & \cA(uY).
}$$
Combining with (i) and claim \ref{cl_split-fibrations}, we deduce
a natural equivalence of fibrations over $\cC$ :
$$
\Fib(v)_*(\cA)\isom\Fib(u)^*(\cA)
$$
and it is easily seen that the system of such equivalences amounts
to a pseudo-natural equivalence of pseudo-functors :
$$
\Fib(v)_*\isom\Fib(u)^*.
$$
\end{remark}

\sset\subsubsection{}\label{subsec_upgrade-to-Fib_!}
For every small category $\cC$ we have a natural functor
$$
\cC/-:\cC\to\Fib(\cC)
\qquad
X\mapsto(\ss_X:\cC/X\to\cC)
$$
that assigns to every morphism $f:X\to Y$ of $\cC$ the
(cartesian) functor $f_*:\cC/X\to\cC/Y$ as in
\eqref{eq_push-for}. The following proposition upgrades
remark \ref{rem_was-cofinal}(ii) from presheaves to
fibrations :

\begin{proposition}\label{prop_upgrade-to-Fib_!}
Let $u:\cC\to\cB$ be any functor between small categories.
With the notation of \eqref{subsec_upgrade-to-Fib_!}, we have
a pseudo-commutative diagram of\/ $2$-categories :
$$
\xymatrix@C+20pt{
\cC \ar[r]^-u \ar[d]_{\cC/-} & \cB \ar[d]^{\cB/-} \\
\Fib(\cC) \ar[r]^-{\Fib(u)_!} & \Fib(\cB).
}$$
\end{proposition}
\begin{proof} Set $\Fib(u)_!(\cC/-):=\Fib(u)_!\circ(\cC/-):
\cC\to\Fib(\cB)$ and let $\cB/u-:\cC\to\Fib(\cB)$ be the
composition of $u$ with the pseudo-funtor $\cB/-$. To ease
notation, we set :
$$
\begin{aligned}
F&\,:=\sCart_\cB(\Fib(u)_!(\cC/-),-):\cC^o\times\Fib(\cB)\to\bCat \\
G&\,:=\sCart_\cB(\cB/u-,-):\cC^o\times\Fib(\cB)\to\bCat
\end{aligned}
$$
(notation of example \ref{ex_first-from-Cats-to-PsFuns}(ii)).
Let $\theta$ be the $2$-adjunction for the pair
$(\Fib(u)_!,\Fib(u)^*)$; we deduce a pseudo-natural equivalence
$$
\omega:=\theta*((\cC/-)^o\times\one_{\Fib(\cB)}):
F\Rightarrow\sCart_\cC(\cC/-,\Fib(u)^*).
$$
On the other hand, we have a strict pseudo-natural equivalence :
$$
\lambda:G\Rightarrow\sCart_\cC(\cC/-,\Fib(u)^*).
$$
Namely, for every $X\in\Ob(\cC)$ and every $\cB$-fibration
$\cA$, the functor
$$
\lambda_{X,\cA}:\cA(uX)\to\Fib(u)^*\cA(X)
$$
assigns to every cartesian section $\phi\in\cA(uX)$ the unique
cartesian section $\phi^*\in\cA\times_\cB\cC(X)$ whose composition
with the projection $\pi:\cA\times_\cB\cC\to\cA$ equals
$\phi\circ u_{|X}$.
To every natural $\cB$-transformation $\alpha:\phi\Rightarrow\phi'$
between such functors, $\lambda_{X,\cA}$ assigns the unique natural
$\cC$-transformation $\alpha^*$ such that
$\pi*\alpha^*=\alpha*u_{|X}$. We may then pick a quasi-inverse
$\mu$ for $\lambda$, and consider the pseudo-natural equivalence
$$
\gamma:=\mu\odot\omega:F\isom G.
$$
For every $X\in\Ob(\cC)$, let $i_X:\Fib(\cB)\to\cC^o\times\Fib(\cB)$
be the unique strict pseudo-functor whose composition with
the projection $\cC^o\times\Fib(\cB)\to\Fib(\cB)$ equals
$\one_{\Fib(\cB)}$, and whose composition with the projection
$\cC^o\times\Fib(\cB)\to\cC^o$ is the constant pseudo-functor
$\sF_{X^o}$ (notation of \eqref{sec_pseudo-cones}); we deduce
the pseudo-natural equivalence
$$
\gamma*i_X:\sCart_\cB(\Fib(u)_!(\cC/X),-)\isom\sCart_\cB(\cB/uX,-).
$$
We show more generally :

\begin{claim}\label{cl_tordu-yoneda}
Let $\cF$ and $\cG$ be two $\cB$-fibrations, and
$\beta_\bullet:\sCart_\cB(\cF,-)\isom\sCart_\cB(\cG,-)$ a
pseudo-natural equivalence. Then $\beta^*:=\beta_\cF(\one_\cF)$
is a $\cB$-equivalence of categories $\cG\isom\cF$.
\end{claim}
\begin{pfclaim} Let
$\alpha_\bullet:\sCart_\cB(\cG,-)\to\sCart_\cB(\cF,-)$ be a
quasi-inverse of $\beta_\bullet$, and set
$\alpha^*:=\alpha_\cG(\one_\cG):\cF\to\cG$. There follow
essentially commutative diagrams :
$$
\xymatrix{
\sCart_\cB(\cF,\cF) \ar[rr]^-{\sCart_\cB(\cF,\alpha^*)}
\ar[d]_{\beta_\cF} & & \sCart_\cB(\cF,\cG) \ar[d]^{\beta_\cG} &
\sCart_\cB(\cG,\cG) \ar[rr]^-{\sCart_\cB(\cG,\beta^*)}
\ar[d]_{\alpha_\cG} & & \sCart_\cB(\cG,\cF) \ar[d]^{\alpha_\cF} \\
\sCart_\cB(\cG,\cF) \ar[rr]^-{\sCart_\cB(\cG,\alpha^*)} & &
\sCart_\cB(\cG,\cG) & \sCart_\cB(\cF,\cG)
\ar[rr]^-{\sCart_\cB(\cF,\beta^*)} & & \sCart_\cB(\cF,\cF).
}$$
From the left diagram we get isomorphisms :
$\alpha^*\circ\beta^*\isom\beta_\cG(\alpha^*)=
\beta_\cG\circ\alpha(\one_\cG)\isom\one_\cG$. Likewise, from
the right diagram we get an isomorphism
$\beta^*\circ\alpha^*\isom\one_\cF$, whence the contention.
\end{pfclaim}

To ease notation, for every $X\in\Ob(\cC)$ and every
morphism $f:X\to Y$ of $\cC$ set
$$
[X]:=\Fib(u)_!(\cC/X)
\qquad\text{and}\qquad
[f]:=\Fib(u)_!(f_*):[X]\to[Y]
$$
According to claim \ref{cl_tordu-yoneda}, for every
$X\in\Ob(\cC)$ we get an equivalence of $\cB$-categories :
$$
\Gamma_X:=\gamma_{X,[X]}(\one_{[X]}):\cB/uX\isom[X].
$$
To conclude, it remains to show that the rule :
$X\mapsto\Gamma_X$ yields a pseudo-natural equivalence
as sought. To this aim, we need to exhibit a coherence
constraint for $\Gamma$. Now, notice that, after possibly
replacing it by a pseudo-naturally isomorphic pseudo-functor,
we may assume that $\Fib(u)_!$ is unital (proposition
\ref{prop_towards-2-yoneda}); then, for every morphism
$f:X\to Y$ in $\cC$ consider the diagram of oriented squares :
$$
\xymatrix@C+40pt{
F(Y,[Y]) \ar[r]^-{\gamma_{Y,[Y]}}
\drtwocell\omit{_\ \ \ \ \ \ \ \ \tau^\gamma_{f,\one_{[Y]}}}
\ar[d]_{F(f,\one_{[Y]})} & G(Y,[Y]) \ar[d]^{G(f,\one_{[Y]})} \\
F(X,[Y]) \ar[r]|-{\gamma_{X,[Y]}} & G(X,[Y]) \\
F(X,[X]) \ar[r]_-{\gamma_{X,[X]}}
\ar[u]^{F(\one_X,[f])} & \ar[u]_{G(\one_X,[f])} G(X,[X]).
\ultwocell\omit{_\qquad\qquad\tau^\gamma_{\one_X,[f]}} \\
}$$
where $\tau^\gamma_{\bullet,\bullet}$ is the coherence constraint of
$\gamma$. We deduce isomorphisms of $\cB$-cartesian functors :
$$
[f]\circ\Gamma_X\xrightarrow{(\tau^\gamma_{\one_X,[f]})_{\one_{[X]}}}
\gamma_{X,[Y]}([f])\xleftarrow{(\tau^\gamma_{f,\one_{[Y]}})_{\one_{[Y]}}}
\Gamma_Y\circ(uf)_*
$$
so our candidate coherence constraint is :
$$
\tau^\Gamma_f:=(\tau^\gamma_{f,\one_{[Y]}})_{\one_{[Y]}}^{-1}
\odot(\tau^\gamma_{\one_X,[f]})_{\one_{[X]}}:
[f]\circ\Gamma_X\isom\Gamma_Y\circ(uf)_*.
$$
With this definition, it is already clear that
$\tau^\Gamma_{\one_X}=\one_{\Gamma_X}$ for every $X\in\Ob(\cC)$.
Next, let $f:X\to Y$ and $g:Y\to Z$ be two morphisms of
$\cC$; we need to check that :
$$
(\tau^\Gamma_g*(uf)_*)\odot([g]*\tau^\Gamma_f)=
\tau^\Gamma_{g\circ f}\odot(\gamma^{\Fib(u)_!}_{f_*,g_*}*\Gamma_X)
$$
where $\gamma^{\Fib(u)_!}$ denotes the coherence constraint of
$\Fib(u)_!$
But the coherence condition for $\tau^\gamma$, relative to
the composition of $1$-cells
$(\one_{X^o},[g])\circ(f^o,\one_{[Y]})=(f^o,[g])=
(f^o,\one_{[Z]})\circ(\one_{Y^o},[g])$ yields the identities :
$$
(\tau^\gamma_{\one_X,[g]})_{[f]}\odot([g]*(\tau^\gamma_{f,\one_{[Y]}}))=
(\tau^\gamma_{f,[g]})_{\one_{[Y]}}=
(\tau^\gamma_{f,\one_{[Z]}})_{[g]}\odot(\tau^\gamma_{\one_Y,[g]}*(uf)_*)
$$
By the same token, from the identity :
$(f^o,\one_{[Z]})\circ(g^o,\one_{[Z]})=((g\circ f)^o,\one_{[Z]})$
we get :
$$
\gamma_{X,[Z]}(\gamma^{\Fib(u)_!}_{f_*,g_*})\odot
(\tau^\gamma_{f,\one_{[Z]}})_{[g]}\odot
((\tau^\gamma_{g,\one_{[Z]}})_{\one_{[Z]}}*(uf)_*)=
(\tau^\gamma_{g\circ f,\one_{[Z]}})_{\one_{[Z]}}
$$
whence :
$$
\begin{aligned}
(\tau^\Gamma_g*(uf)_*)\!\odot\!([g]*\tau^\Gamma_f)&\!=\!
((\tau^\gamma_{g,\one_{[Z]}})^{-1}_{\one_{[Z]}}*(uf)_*)\!\odot\!
(\tau^\gamma_{f,\one_{[Z]}})_{[g]}^{-1}\!\odot\!(\tau^\gamma_{\one_X,[g]})_{[f]}
\!\odot\!([g]\!*\!(\tau^\gamma_{\one_X,[f]})_{\one_{[X]}}) \\
&\!=\!(\tau^\gamma_{g\circ f,\one_{[Z]}})_{\one_{[Z]}}^{-1}\odot
\gamma_{X,[Z]}(\gamma^{\Fib(u)_!}_{f_*,g_*})\!\odot\!
(\tau^\gamma_{\one_X,[g]})_{[f]}\!\odot\!
([g]\!*\!(\tau^\gamma_{\one_X,[f]})_{\one_{[X]}}).
\end{aligned}
$$
Thus, we are reduced to checking that :
$$
(\tau^\gamma_{\one_X,[g\circ f]})_{\one_{[X]}}\odot
(\gamma^{\Fib(u)_!}_{f_*,g_*}*\Gamma_X)=
\gamma_{X,[Z]}(\gamma^{\Fib(u)_!}_{f_*,g_*})\odot
(\tau^\gamma_{\one_X,[g]})_{[f]}\odot
([g]*(\tau^\gamma_{\one_X,[f]})_{\one_{[X]}}).
$$
But the latter follows once again from the coherence condition
for $\tau^\gamma$, applied to the identity :
$(\one_X,[g])\circ(\one_X,[f])=(\one_X,[g\circ f])$.
\end{proof}

In the same vein, we point out the following :

\begin{proposition}\label{prop_same-vein-Fib}
Let $\cC$, $\cC'$ be two small categories, $u:\cC\to\cC'$
a functor, and $\sV$ a universe such that $\cC^\wedge$ and
$\cC'^\wedge$ are $\sV$-small. We have a pseudo-commutative
diagram of\/ $2$-categories :
$$
\xymatrix@C+40pt{
\sV\tdu\Fib(\cC'^\wedge_\sU) \ar[r]^-{\sV\tdu\Fib(u^\wedge_\sU)_!}
\ar[d]_{\sV\tdu\Fib(h_{\cC'})^*} &
\sV\tdu\Fib(\cC^\wedge_\sU) \ar[d]^{\sV\tdu\Fib(h_\cC)^*} \\
\sV\tdu\Fib(\cC') \ar[r]^-{\sV\tdu\Fib(u)^*} & \sV\tdu\Fib(\cC).
}$$
\end{proposition}
\begin{proof} Let $\cE'$ be any fibration with $\sV$-small
fibres over $\cC'^\wedge$; pick a unital cleavage for $\cE'$
and let $\sc'$ be its associated pseudo-functor. Set
$\cE:=\sV\tdu\Fib(u^\wedge_\sU)_!\cE'$; by inspecting
\eqref{subsec_adjoints-of-Fib(rho)} we see that for every
$F\in\Ob(\cC^\wedge)$ the fibre category $\cE_F$ represents
the strong $2$-colimit of
$$
\sc'\circ\st^o_F:(F/u^\wedge\cC'^\wedge)\to\sV\tdu\bCat
$$
where $\st_F:(F/u^\wedge\cC'^\wedge)\to\cC'^\wedge$ denotes the
target functor. Fix a universal pseudo-cocone
$$
e^F_\bullet:\sc'\circ\st^o_F\Rightarrow\sF_{\cE_F}
\qquad
(F\xrightarrow{\phi}u^\wedge G)\mapsto
(\cE'_G\xrightarrow{e^F_\phi}\cE_F)
$$
and let $\tau^F_\bullet$ be the coherence constraint of
$e^F_\bullet$. Then $\cE$ admits a split cleavage $\sd$,
that assigns to every morphism $\psi:F'\to F$ of $\cC^\wedge$
the unique functor $\cE_\psi:\cE_F\to\cE_{F'}$ such that
$$
\sF_{\cE_\psi}\odot e^F=e^{F'}*\psi^{*o}
$$
where $\psi^*:F/u^\wedge\cC'^\wedge\to F'/u^\wedge\cC'^\wedge$
is defined as in \eqref{eq_push-for}. Hence,
$\sV\tdu\Fib(h_\cC)^*$ admits the split cleavage
$\sd\circ h^o_\cC$. For every $X\in\Ob(\cC)$, let
$\eta_X:h_X\to u^\wedge h_{uX}$ be the morphism given
by the rule : $\phi\mapsto u(\phi)\in h_{uX}(uY)$ for every
$Y\in\Ob(\cC)$ and every $\phi\in h_X(Y)$. We notice :

\begin{claim}\label{cl_eta-is-initial}
For every $X\in\Ob(\cC)$, the morphism $\eta_\cE$ is an
initial object of $h_X/u^\wedge\cC'^\wedge$.
\end{claim}
\begin{pfclaim} Indeed,  In light of
remark \ref{rem_was-cofinal}(ii,iii), we may regard
$\eta_X$ as the unit of the natural adjunction for
the adjoint pair $(u_!,u^\wedge)$. Hence, for every
$G\in\Ob(\cC'^\wedge)$ and every morphism of presheaves
$\beta:h_X\to u^\wedge G$ there exists a unique morphism
of presheaves $\beta':h_{uX}\to G$ on $\cC'$ such that
$\beta=u^\wedge(\beta')\circ\eta_X$ (explicitly, under
the canonical bijection provided by Yoneda's lemma,
$\beta$ corresponds to a unique section $s_\beta\in(u^\wedge G)X$,
and $\beta'$ corresponds then to the same $s_\beta$,
regarded as an element of $G(uX)$). The claim is an
immediate consequence.
\end{pfclaim}

From claim \ref{cl_eta-is-initial} and proposition
\ref{prop_pseudo-initial}, we get an equivalence of categories
$$
\eps_X:=e^{h_X}_{\eta_X}:\cE'_{h_{uX}}\to\cE_{h_X}
\qquad
\text{for every $X\in\Ob(\cC)$}.
$$

\begin{claim}\label{cl_long-knee}
The rule : $X\mapsto\eps_X$ yields a pseudo-natural equivalence
$\eps:\sc'\circ(h_{\cC'}\circ u)^o\isom\sd\circ h^o_\cC$.
\end{claim}
\begin{pfclaim} For every morphism $f:Y\to X$ in $\cC$ we have a
natural isomorphism :
$$
\tau^\eps_f:\sd_{h_f}\circ\eps_X=e^{h_Y}_{\eta_X\circ h_f}=
e^{h_Y}_{u^\wedge(h_{uf})\circ\eta_Y}
\xrightarrow{\tau^{h_Y}_{h_{uf}}}\eps_Y\circ\sc'_{h_{uf}}
$$
and we check that the system $(\tau^\eps_f~|~f\in\rMorph(\cC))$
yields a coherence constraint for $\eps$. First, we get
$\tau^\eps_{\one_X}=\one_{\eps_X}$ for every $X\in\Ob(\cC)$,
by remark \ref{rem_unital}(ii). Next, notice that for every
morphism $\psi:F'\to F$ in $\cC^\wedge$ we have :
$$
\cE_\psi*\tau^F_\mu=\tau^{F'}_\mu
\qquad
\text{for every morphism
$F'/\mu:(F'\to u^\wedge G_1)\to(F'\to u^\wedge G_2)$
in $F'/u^\wedge\cC'^\wedge$}.
$$
Thus, let $f:X\to Y$ and $g:Y\to Z$ be two morphisms in $\cC$;
we deduce :
$$
\begin{aligned}
(\eps_Z*\gamma^{\sc'}_{h_{uf},h_{ug}})\odot
(\tau^\eps_g*\sc'_{uf})\odot(\sd_{h_g}*\tau^\eps_f)&\,=
(\eps_Z*\gamma^{\sc'}_{h_{uf},h_{ug}})\odot
(\tau^{h_Z}_{h_{ug}}*\sc'_{uf})\odot\tau^{h_Z}_{h_{uf}} \\
&\,=(\eps_Z*\gamma^{\sc'}_{h_{uf},h_{ug}})\odot\tau^{h_Z}_{h_{u(g\circ f)}} \\
&\,=\tau^\eps_{g\circ f}
\end{aligned}
$$
whence the contention.
\end{pfclaim}

From claim \ref{cl_long-knee} we deduce an equivalence of
fibrations over $\cC$ :
$$
\Omega_{\cE'}:\sV\tdu\Fib(h_{\cC'}\circ u)^*\cE'\isom
\sV\tdu\Fib(h_\cC)^*\circ\sV\tdu\Fib(u^\wedge_\sU)_!\cE'.
$$
Next, let $\phi':\cE'_1\to\cE'_2$ be a cartesian functor
of fibrations over $\cC'^\wedge$; let
$\cE_i:=\sV\tdu\Fib(u^\wedge)_!\cE'_i$ for $i=1,2$, and
set $\phi:=\sV\tdu\Fib(u^\wedge_\sU)_!\phi':\cE_1\to\cE_2$.
For $i=1,2$, pick also unital a cleavage $\sc'_i$ for
$\cE'$, so that $\phi'$ corresponds to a pseudo-natural
transformation $\omega':\sc'_1\Rightarrow\sc'_2$. For every
$F\in\Ob(\cC^\wedge)$ we pick universal cocones
$e^F_{1,\bullet}:\sc'_i\circ\st^o_F\Rightarrow\sF_{\cE_{i,F}}$,
and denote by $\sd_1$ and $\sd_2$ the natural split cleavages
for $\cE_1$ and $\cE_2$ described in the foregoing. Then
$\phi$ corresponds to the strict pseudo-natural transformation
$\omega:\sd_1\Rightarrow\sd_2$ that assigns to every
$F\in\Ob(\cC^\wedge)$ the unique functor
$\omega_F:\cE_{1,F}\to\cE_{2,F}$ such that
\set\begin{equation}\label{eq_knee-better}
e^F_{2,\bullet}\odot(\omega'*\st^o_F)=
\sF_{\omega_F}\odot e^F_{1,\bullet}.
\end{equation}
Let $\eps_i:\sc'_i\circ(h_{\cC'}\circ u)^o\isom\sd_i\circ h^o_\cC$
be the pseudo-natural equivalence of claim \ref{cl_long-knee},
for $i=1,2$.

\begin{claim}\label{cl_knee-better}
$(\omega*h^o_\cC)\odot\eps_1=\eps_2\odot(\omega'*(h_{\cC'}\circ u)^o)$.
\end{claim}
\begin{pfclaim} Directly from \eqref{eq_knee-better} we see that
the two sides of the identity of the claim agree on every
$X\in\Ob(\cC)$. It remains to check that the respective coherence
constraints agree as well. However, for every $F\in\Ob(\cC^\wedge)$
denote by $\tau^F_{i,\bullet}$ the coherence constraint of
$e^F_{i,\bullet}$, for $i=1,2$; the coherence constraint of
$(\omega*h^o_\cC)\odot\eps_1$ assigns to every morphism $f:X\to Y$
of $\cC$ the natural isomorphism of functors
$\omega_{h_X}*\tau^{h_Y}_{1,h_{uf}}$, whereas the coherence constraint
of $\eps_2\odot(\omega'*(h_{\cC'}\circ u)^o)$ assigns to $f$ the
natural isomorphism $\tau^{h_Y}_{2,h_{uf}}*\omega'_{h_{uY}}$. Then
again the required identity follows directly from
\eqref{eq_knee-better}.
\end{pfclaim}

From claim \ref{cl_knee-better} we deduce a commutative diagram
of cartesian functors :
$$
\xymatrix@C+20pt{ \sV\tdu\Fib(h_{\cC'}\circ u)^*\cE'_1
\ar[d]_{\sV\tdu\Fib(h_{\cC'}\circ u)^*\phi'}
\ar[r]^-{\Omega_{\cE'_1}} & \sV\tdu\Fib(h_\cC)^*\cE_1
\ar[d]^{\sV\tdu\Fib(h_\cC)^*\phi} \\
\sV\tdu\Fib(h_{\cC'}\circ u)^*\cE'_2
\ar[r]^-{\Omega_{\cE'_2}} & \sV\tdu\Fib(h_\cC)^*\cE_2.
}$$
Lastly, for $\cE'_1,\cE'_2$ as in the foregoing, let
$\phi'_1,\phi'_2:\cE'\to\cE'_2$ be two cartesian functors,
$\beta':\phi'_1\Rightarrow\phi'_2$ a natural
$\cC'^\wedge$-transformation, and set
$\phi_i:=\sV\tdu\Fib(u^\wedge_\sU)_!\phi'_i$ for $i=1,2$ and
$\beta:=\sV\tdu\Fib(u^\wedge_\sU)_!\beta':\phi_1\Rightarrow\phi_2$.
Say that $\phi'_1$ and $\phi'_2$ correspond to pseudo-natural
transformations $\omega'_1,\omega'_2:\sc'_1\Rightarrow\sc'_2$,
and $\beta'$ corresponds to a modification
$\Xi':\omega'_1\leadsto\omega'_2$; then $\phi_1$ and $\phi_2$
correspond to strict pseudo-natural transformations
$\omega_1,\omega_2:\sd_1\Rightarrow\sd_2$ described as in the
foregoing, and $\beta$ corresponds to the modification
$\Xi:\omega_1\leadsto\omega_2$ assigning to every
$F\in\Ob(\cC^\wedge)$ the unique natural transformation
$\Xi_F:\omega_{1,F}\Rightarrow\omega_{2,F}$ with :
$$
e^F_{2,\bullet}*(\Xi'\circ\st^o_F)=\sF_{\Xi_F}*e^F_{1,\bullet}
$$
which specializes to the identity :
$\Xi_{h_X}\odot e^{h_X}_{1,\eta_X}=e^{h_X}_{2,\eta_X}\odot\Xi'_{h_{uX}}$
for every $X\in\Ob(\cC)$. The latter in turns yields the identity :
$$
(\Xi\circ h^o_\cC)*\eps_1=\eps_2*(\Xi'\circ(h_{\cC'}\circ u)^o)
$$
which finally shows that :
$$
\sV\tdu\Fib(h_\cC)^*(\beta)*\Omega_{\cE'_1}=
\Omega_{\cE'_2}*\sV\tdu\Fib(h_{\cC'}\circ u)^*(\beta').
$$
Summing up, this shows that the system of equivalences
$\Omega_\bullet$ yields a strict pseudo-natural equivalence
of strict pseudo-functors $\sV\tdu\Fib(h_{\cC'}\circ u)^*\isom
\sV\tdu\Fib(h_\cC)^*\circ\sV\tdu\Fib(u^\wedge_\sU)_!$, as stated.
\end{proof}

\sset\subsubsection{}\label{subsec_lims-in-Fib}
We conclude this section by showing that the $2$-limits and
$2$-colimits are computed fibrewise in the $2$-category
$\Fib(\cB)$, a result that extends the corresponding assertion
for the category of presheaves (corollary \ref{cor_pre-misc}(ii)).
To begin with, let $\cB$ be any category; with every $B\in\Ob(\cB)$
we associate a strict pseudo-functor
$$
\fib_B:\Fib(\cB)\to\bCat
\qquad
\cA\mapsto\cA_B.
$$
If $\cA$ and $\cA'$ are two $\cB$-fibrations,
$\psi,\psi':\cA\to\cA'$ two $\cB$-cartesian functors, and
$\beta:\psi\Rightarrow\psi'$ a natural $\cB$-transformation,
then $\fib_B(\psi):\cA_B\to\cA'_B$ is the restriction
of $\psi$, and
$\fib_B(\beta):\fib_B(\psi)\Rightarrow\fib_B(\psi')$ is the
restriction of $\beta$. We may then state :

\begin{theorem}\label{th_fib-limits-are-fibrewise}
With the notation of \eqref{subsec_lims-in-Fib}, we have :
\begin{enumerate}
\item
The $2$-category $\Fib(\cB)$ is strongly $2$-complete and
strongly $2$-cocomplete.
\item
For every $B\in\Ob(\cB)$ the pseudo-functor $\fib_B$ commutes
with the strong $2$-limit and strong $2$-colimit of every
pseudo-functor $I\to\Fib(\cB)$ from any small $2$-category $I$. 
\end{enumerate} 
\end{theorem}
\begin{proof} We prove the $2$-completeness of $\Fib(\cB)$; the
$2$-cocompleteness follows by the same argument, considering
the opposite $2$-categories : the details shall be left to the
reader. In view of theorem \ref{th_fundamental-fibrations},
it suffices to show the corresponding assertions for the
$2$-category $\sPsFun(\cB^o,\bCat)$ and the similar strict
pseudo-functors that we denote as well by
$$
\fib_B:\sPsFun(\cB^o,\bCat)\to\bCat
\qquad
\sc\mapsto\sc_B.
$$
Notice that every morphism $g:B'\to B$ of $\cB$ induces
a pseudo-natural transformation
$$
\fib_g:\fib_B\Rightarrow\fib_{B'}
\qquad
\sc\mapsto(\sc_g:\sc_B\to\sc_{B'})
$$
with coherence constraint given by the rule :
$$
\tau^{\fib_g}_\beta:=(\tau^\beta_g)^{-1}
$$
for every pair of pseudo-natural functors $\sc,\sd:\cB^o\to\bCat$
and every pseudo-natural transformation $\beta:\sc\Rightarrow\sd$
with coherence constraint $\tau^\beta$. Moreover, if $h:B''\to B'$
is any other morphism of $\cB$, we get the invertible modification
$$
\Gamma_{h,g}:\fib_h\odot\fib_g\leadsto\fib_{g\circ h}
\qquad
\sc\mapsto\gamma^\sc_{h,g}
$$
where $(\delta^\sc,\gamma^\sc)$ denotes the coherence constraint of
$\sc$. Indeed, it is easily seen that the compatibility condition
for $(\Gamma_{h,g})_\beta$ corresponding to a pseudo-natural
transformation $\beta:\sc\Rightarrow\sd$ is equivalent to
the coherence axiom for $\tau^\beta$ : details left to the reader.
Likewise, for every third morphism $k:B'''\to B''$ of $\cB$, the
composition axiom for pseudo-functors translates as :
\set\begin{equation}\label{eq_translates-comp-axiom}
\Gamma_{k,g\circ h}\odot(\fib_k*\Gamma_{h,g})=
\Gamma_{h\circ k,g}\odot(\Gamma_{k,h}*\fib_g).
\end{equation}
Let now $F:I\to\sPsFun(\cB^o,\bCat)$ be any pseudo-functor;
we need to exhibit a $2$-limit for $F$, and by virtue of
proposition \ref{prop_towards-2-yoneda}, we may assume that
$Fi:\cB^o\to\bCat$ is unital for every $i\in\Ob(I)$. Then, for
every $B\in\Ob(\cB)$ choose a strong $2$-limit $(L(B),\pi^B)$
of the pseudo-functor $\fib_B\circ F$ (theorem
\ref{th_bCat-cpt}), and let $\tau^B$ be the coherence
constraint of the pseudo-cone
$\pi^B:\sF_{L(B)}\Rightarrow\fib_B\circ F$. Let $g:B'\to B$
be any morphism of $\cB$; by the (strong) universality of
$\pi^{B'}$, we find a unique functor $L(g):L(B)\to L(B')$
such that
\set\begin{equation}\label{eq_use-strong-lims}
\pi^{B'}\odot\sF_{L(g)}=(\fib_g*F)\odot\pi^B.
\end{equation}
Notice that $\fib_{\one_B}*F=\one_{\fib_B\circ F}$, since $Fi$
is unital for every $i\in\Ob(I)$; therefore :
$$
L(\one_B)=\one_{L(B)}
\qquad
\text{for every $B\in\Ob(\cB)$}.
$$
Moreover, if $h:B''\to B'$ is another morphism of $\cB$,
we have identities :
$$
\begin{aligned}
\pi^{B''}\odot\sF_{L(h)\circ L(g)}&\,=
(\fib_h*F)\odot(\fib_g*F)\odot\pi^B=
((\fib_h\odot\fib_g)*F)\odot\pi^B \\
\pi^{B''}\odot\sF_{L(h\circ g)}&\,=(\fib_{h\circ g})\odot\pi^B
\end{aligned}
$$
so there exists a unique natural transformation
$\gamma^L_{h,g}:L(h)\circ L(g)\Rightarrow L(g\circ h)$
such that
$$
\pi^{B''}*\sF_{\gamma^L_{h,g}}=(\Gamma_{h,g}\circ F)*\pi^B.
$$
\begin{claim} The rules $B\mapsto L(B)$ and
$(g:B'\to B)\mapsto L(g)$ for every $B\in\Ob(\cB)$ and
every morphism $g$ of $\cB$ define a unital pseudo-functor
$L:\cB^o\to\bCat$ with coherence constraint given by
the system of natural transformations
$\gamma^L_{\bullet,\bullet}$.
\end{claim}
\begin{pfclaim} The unit axiom for $\gamma^L$ is clear
from the construction. To check the composition axiom,
let $g$ and $h$ be as in the foregoing, and consider as
well a third morphism $k:B'''\to B''$ of $\cB$. Since
$\pi^{B'''}$ is a universal pseudo-cone, it suffices to
show that
$$
X:=\pi^{B'''}*(\sF_{\gamma^L_{k,g\circ h}}\odot\sF_{L(k)*\gamma^L_{h,g}})=
Y:=\pi^{B'''}*(\sF_{\gamma^L_{h\circ k,g}}\odot\sF_{\gamma^L_{k,h}*L(g)}).
$$
However, by unwinding the definitions we find :
$$
\begin{aligned}
X=&\,((\Gamma_{k,g\circ h}\circ F)*\pi^B)
\odot(\pi^{B'''}*\sF_{L(k)}*\sF_{\gamma^L_{h,g}}) \\
=&\,((\Gamma_{k,g\circ h}\circ F)*\pi^B)\odot
((\fib_k*F)*\pi^{B''}*\sF_{\gamma^L_{h,g}}) \\
=&\,((\Gamma_{k,g\circ h}\circ F)*\pi^B)\odot
((\fib_k*F)*(\Gamma_{h,g}\circ F)*\pi^B) \\
=&\,((\Gamma_{k,g\circ h}\circ F)*\pi^B)\odot
(((\fib_k*\Gamma_{h,g})\circ F)*\pi^B) \\
Y=&\,
((\Gamma_{h\circ k,g}\circ F)*\pi^B)\odot
(\pi^{B'''}*\sF_{\gamma^L_{k,h}}*\sF_{L(g)}) \\
=&\,((\Gamma_{h\circ k,g}\circ F)*\pi^B)\odot
((\Gamma_{k,h}\circ F)*\pi^{B'}*\sF_{L(g)}) \\
=&\,((\Gamma_{h\circ k,g}\circ F)*\pi^B)\odot
((\Gamma_{k,h}\circ F)*(\fib_g*F)*\pi^B) \\
=&\,((\Gamma_{h\circ k,g}\circ F)*\pi^B)\odot
((\Gamma_{k,h}*\fib_g)\circ F)\odot\pi^B
\end{aligned}
$$
so it suffices to show:
$$
(\Gamma_{\!h\circ k,g}\circ F)\odot((\Gamma_{k,h}*\fib_g)\circ F)=
(\Gamma_{k,g\circ h}\circ F)\odot((\fib_k*\Gamma_{h,g})\circ F)
$$
which follows straightforwardly from
\eqref{eq_translates-comp-axiom}.
\end{pfclaim}

Next, we notice that for every $i\in\Ob(I)$ the rule :
$B\mapsto(\pi^B_i:L(B)\to Fi(B))$ defines a strict
pseudo-natural transformation
$$
\pi_i:L\Rightarrow Fi.
$$
Indeed, the first coherence axiom for $\pi_i$ is easily
checked, recalling that both $L$ and $Fi$ are unital,
and the second one follows directly from the definition
of $\gamma^L$ : details left to the reader.

\begin{claim}\label{cl_switch}
(i)\ \
For every $1$-cell $\phi:i\to j$ in $I$, the rule :
$$
B\mapsto\tau^B_\phi
$$
defines an invertible modification
$\Xi^\phi:F(\phi)\odot\pi_i\leadsto\pi_j$.

(ii)\ \
The system $(\pi_i~|~i\in\Ob(I))$ defines a pseudo-cone
$$
\pi:\sF_L\Rightarrow F
$$
with coherence constraint given by the invertible
modifications $\Xi^\bullet$.
\end{claim}
\begin{pfclaim}(i): Let $\tau^{F\phi}$ and $\tau^{\fib_g*F}$ be the
coherence constraints of $F\phi$ and respectively $\fib_g*F$.
We need to verify the compatibility condition :
$$
Fj(g)*\tau^B_\phi=(\tau^{B'}_j*L(g))\odot(\tau^{F\phi}_g*\pi^B_i)
$$
for every morphism $g:B'\to B$ of $\cB$. However, comparing
the coherence constraints of the two sides of
\eqref{eq_use-strong-lims} we get :
$$
(Fj(g)*\tau^B_\phi)\odot(\tau^{\fib_g*F}_\phi*\pi^B_i)=
\tau^{B'}_\phi*L(g).
$$
Since $\tau^{\fib_g*F}_\phi=(\tau^{F\phi}_g)^{-1}$, the assertion
follows.

(ii): Denote by $(\delta^F,\gamma^F)$ the coherence constraint
of $F$; we need to check the identities
$$
\Xi^{\one_i}\odot(\delta^F_i*\pi_i)=\one_{\pi_i}
\qquad
\Xi^\psi\odot(F\psi*\Xi^\phi)=
\Xi^{\psi\circ\phi}\odot(\gamma^F_{\phi,\psi}*\pi_i)
$$
for every $i,j,k\in\Ob(I)$ and every pair of $1$-cells
$\phi:i\to j$ and $\psi:j\to k$ of $I$. However, the
latter are none else than a rephrasing of the coherence
axioms for the pseudo-cones $\pi^B$, for $B$ ranging over
all the objects of $\cB$.
\end{pfclaim}

It remains to check that the pseudo-cone
$\pi:\sF_L\Rightarrow F$ of claim \ref{cl_switch}(ii)
is universal. To this aim, let $X:\cB^o\to\bCat$ be
any pseudo-functor, and $\beta:\sF_X\Rightarrow F$
any pseudo-cone; we need to exhibit a pseudo-natural
trasformation $t:X\Rightarrow  L$ such that
$\pi\odot\sF_t=\beta$. Let $\omega:X\isom X^u$ be the
isomorphism furnished by proposition
\ref{prop_towards-2-yoneda}, with $X^u$ a unital
pseudo-functor; we get a commutative diagram of categories :
$$
\xymatrix{ \sPsNat(X^u,L) \ar[r] \ar[d]_{\sPsNat(\omega,L)} &
\sPsNat(\sF_{X^u},F) \ar[d]^{\sPsNat(\sF_\omega,F)} \\
\sPsNat(X,L) \ar[r] & \sPsNat(\sF_X,F)
}$$
whose vertical arrows are isomorphisms. We need to check
that the bottom horizontal arrow is an isomorphism, so it
suffices to show that the same holds for the top horizontal
arrow; we may thus assume that $X$ is unital.
Now, $\fib_B*\beta$ is a pseudo-cone with vertex $X(B)$
and basis $\fib_B\circ F:I\to\bCat$; then there exists
a unique functor $t_B:X(B)\to L(B)$ such that
$$
\pi^B\odot\sF_{t_B}=\fib_B*\beta
\qquad
\text{for every $B\in\Ob(\cB)$}.
$$
Denote by $\tau^\beta$ the coherence constraint of $\beta$,
and by $\tau^{\beta_i}$ the coherence constraint of the
pseudo-natural transformation $\beta_i:X\Rightarrow Fi$,
for every $i\in\Ob(I)$; we remark :

\begin{claim}\label{cl_fresh-books}
For every morphism $g:B'\to B$ of $\cB$, the rule :
$$
i\mapsto\tau^{\beta_i}_g
\qquad
\text{for every $i\in\Ob(I)$}
$$
defines an invertible modification
$\Lambda^g:(\fib_g*F)\odot(\fib_B*\beta)\leadsto
(\fib_{B'}*\beta)\odot\sF_{X(g)}$.
\end{claim}
\begin{pfclaim} The coherence constraint of
$(\fib_g*F)\odot(\fib_B*\beta)$ assigns to every $1$-cell
$\phi:i\to j$ of $I$ the composition of oriented squares :
$$
\xymatrix{ X(B) \ddouble \ar[rr]^-{\beta_{i,B}}
\drrtwocell\omit{_\ \ \ \ \ \tau^\beta_{\phi,B}} & &
\drrtwocell\omit{_\ \ \ \ \ \ \ \ \ (\tau^{F\phi}_g)^{-1}}
Fi(B) \ar[rr]^-{Fi(g)} \ar[d]|{F\phi(B)} & &
Fi(B') \ar[d]^{F\phi(B')} \\
X(B) \ar[rr]_-{\beta_{j,B}} & & Fj(B) \ar[rr]_-{Fj(g)} & &
Fj(B')
}$$
whereas the coherence constraint of
$(\fib_{B'}*\beta)\odot\sF_{X(g)}$ assigns to $\phi$ the
natural transformation
$$
\tau^\beta_{\phi,B'}*X(g):F\phi(B')\circ\beta_{i,B'}\circ X(g)
\Rightarrow\beta_{j,B'}\circ X(g).
$$
Thus, for every such $\phi$ we need to show the identity :
$$
\tau^{\beta_j}_g\odot(Fj(g)*\tau^\beta_{\phi,B})\odot
((\tau^{F\phi}_g)^{-1}*\beta_{i,B})=
\tau^\beta_{\phi,B'}\odot(F\phi(B')*\tau^{\beta_i}_g).
$$
But the latter is equivalent to the compatibility condition
for the invertible modification
$\tau^\beta_\phi:F\phi\odot\beta_i\leadsto\beta_j$, corresponding
to the morphism $g$.
\end{pfclaim}

Notice now that
$(\fib_g*F)\odot(\fib_B*\beta)=\pi^{B'}\odot\sF_{L(g)\circ t_B}$
and $(\fib_{B'}*\beta)\odot\sF_{X(g)}=\pi^{B'}\odot\sF_{t_{B'}\circ X(g)}$
for every morphism $g:B'\to B$ of $\cB$. In view of claim
\ref{cl_fresh-books}, and since $\pi^{B'}$ is universal, there
exists therefore a unique natural isomorphism of functors
$$
\tau^t_g:L(g)\circ t_B\Rightarrow t_{B'}\circ X(g)
\qquad\text{such that}\qquad
\pi^{B'}*\sF_{\tau^t_g}=\Lambda^g.
$$
Notice that, by virtue of remark \ref{rem_unital}(ii), we
have $\Lambda^{\one_B}=\one_{\fib_B*\beta}$ for every $B\in\Ob(\cB)$,
and therefore $\tau^t_{\one_B}=\one_{t_B}$ for every such $B$.
Let us check that the system $(t_B~|~B\in\Ob(\cB))$ yields
the sought pseudo-natural transformation, with coherence
constraint given by the system of isomorphisms of functors
$\tau^t_\bullet$. Indeed, the first coherence axiom follows
easily from the foregoing; it remains therefore only to show
that for every pair of morphisms $g:B'\to B$ and $h:B''\to B'$
of $\cB$ we have :
$$
(t_{B''}*\gamma^X_{h,g})\odot(\tau^t_h*X(g))\odot(L(h)*\tau^t_g)=
\tau^t_{g\circ h}\odot(\gamma^L_{g,h}*t_B)
$$
where $\gamma^X$ denotes the coherence constraint of $X$.
Using the universality of $\pi^{B''}$, we are then reduced
to showing the identity :
$$
U:=(\pi^{B''}*\sF_{t_{B''}}*\sF_{\gamma^X_{h,g}})\odot
(\Lambda^h*\sF_{X(g)})\odot(\pi^{B''}*\sF_{L(h)}*\sF_{\tau^t_g})=
V:=\Lambda^{g\circ h}\odot(\pi^{B''}*\sF_{\gamma^L_{g,h}}*\sF_{t_B}).
$$
However, we have :
$$
\begin{aligned}
U=&\,(\pi^{B''}*\sF_{t_{B''}}*\sF_{\gamma^X_{h,g}})\odot
(\Lambda^h*\sF_{X(g)})\odot((\fib_h*F)*\pi^{B'}*\sF_{\tau^t_g}) \\
=&\,((\fib_{B''}*\beta)*\sF_{\gamma^X_{h,g}})\odot
(\Lambda^h*\sF_{X(g)})\odot((\fib_h*F)*\pi^{B'}*\sF_{\tau^t_g}) \\
=&\,((\fib_{B''}*\beta)*\sF_{\gamma^X_{h,g}})\odot
(\Lambda^h*\sF_{X(g)})\odot((\fib_h*F)*\Lambda^g) \\
V=&\,\Lambda^{g\circ h}\odot((\Gamma_{h,g}\circ F)*\pi^B*\sF_{t_B}).
\end{aligned}
$$
Thus, it suffices to show :
$$
((\fib_{B''}*\beta)*\sF_{\gamma^X_{h,g}})\odot
(\Lambda^h*\sF_{X(g)})\odot((\fib_h*F)*\Lambda^g)=
\Lambda^{g\circ h}\odot((\Gamma_{h,g}\circ F)*(\fib_B*\beta)).
$$
But a simple inspection shows that the latter identity
translates the coherence axiom for $\beta_i$, relative
to the pair of $1$-cells $(h,g)$, and for $i$ ranging
over all the objects of $I$.

\begin{claim}\label{cl_wait-one-week}
(i)\ \
$\pi_i\odot t=\beta_i$ for every $i\in\Ob(I)$.

(ii)\ \
$\pi\odot\sF_t=\beta$.
\end{claim}
\begin{pfclaim} We need to check that the coherence constraint
$\tau^{\pi_i\odot t}$ of $\pi_i\odot t$ equals $\tau^{\beta_i}$.
But for every morphism $g:B'\to B$ of $\cB$ we have
$\tau^{\pi_i\odot t}_g=\pi^{B'}_i*\tau^t_g=\Lambda^g_i$,
whence the contention.

(ii): In view of (i), it suffices to check that $\beta$
and $\pi\odot\sF_t$ have the same coherence constraints.
But the coherence constraint of $\pi\odot\sF_t$ assigns
to every $1$-cell $\phi:i\to j$ in $I$ the invertible
modification $\Xi^\phi*t$, so the assertion follows by
a direct inspection of the definitions.
\end{pfclaim}

Claim \ref{cl_wait-one-week}(ii) shows that the functor
\set\begin{equation}\label{eq_back-belt}
\sPsNat(X,L)\to\sPsNat(\sF_X,F)
\qquad
t\mapsto\pi\odot\sF_t
\end{equation}
is surjective on objects. To check injectivity on objects,
consider two pseudo-natural transformations $t,t':X\Rightarrow L$
such that $\pi\odot\sF_t=\pi\odot\sF_{t'}$. Then we have
$$
\pi^B\odot\sF_{t_B}=\fib_B*(\pi\odot\sF_t)=
\fib_B*(\pi\odot\sF_{t'})=\pi^B\odot\sF_{t'_B}
\qquad
\text{for every $B\in\Ob(\cB)$}.
$$
By strong universality of $\pi^B$, we deduce that $t_B=t'_B$
for every such $B$. Moreover, recall that we have attached
to $\pi\odot\sF_t$ and every morphism $g:B'\to B$ an invertible
modification $\Lambda^g:\pi^{B'}\odot\sF_{L(g)\circ t_B}\leadsto
\pi^{B'}\odot\sF_{t_{B'}\circ X(g)}$. Let $\tau^t$ and $\tau^{t'}$
be the coherence constraint of $t$ and $t'$; by inspecting
the definition, we find that
$$
\Lambda^g_i=(\tau^B_g*t_B)\odot(\pi^{B'}_i*\tau^t_g)
\qquad
\text{for every $i\in\Ob(I)$}.
$$
Since $\pi\odot\sF_t=\pi\odot\sF_{t'}$ and $t_B=t'_B$, we
get $\pi^{B'}_i*\tau^t_g=\pi^{B'}_i*\tau^{t'}_g$ for every
$i\in\Ob(I)$, {\em i.e.}
$\pi^{B'}*\sF_{\tau^t_g}=\pi^{B'}*\sF_{\tau^{t'}_g}$, whence
$\tau^t_g=\tau^{t'}_g$, by the universality of $\pi^{B'}$.
Thus, $t=t'$ as required.

To check that \eqref{eq_back-belt} is faithful, consider
pseudo-natural transformations $t,t':X\Rightarrow L$ and
modifications $\Delta,\Delta':t\leadsto t'$ with
$\pi*\sF_\Delta=\pi*\sF_{\Delta'}$; we deduce:
$$
\pi^B*(\fib_B\circ\sF_\Delta)=\fib_B\circ(\pi*\sF_\Delta)=
\fib_B\circ(\pi*\sF_{\Delta'})=\pi^B*(\fib_B\circ\sF_{\Delta'})
\qquad
\text{for every $B\in\Ob(\cB)$}
$$
and clearly $\fib_B\circ\sF_\Delta$ and $\fib_B\circ\sF_{\Delta'}$
are the constant modifications $\sF_{t_B}\leadsto\sF_{t'_B}$
associated with $\Delta$ and $\Delta'$. By the universality
of $\pi^B$, it follows that
$\fib_B\circ\sF_\Delta=\fib_B\circ\sF_{\Delta'}$ for every
$B\in\Ob(\cB)$, whence $\Delta=\Delta'$, as required.

Lastly, let $\Delta:\pi\odot\sF_t\leadsto\pi\odot\sF_{t'}$
be any modification; for every $B\in\Ob(\cB)$ we obtain
the modification
$\fib_B\circ\Delta:\pi^B\odot\sF_{t_B}\leadsto\pi^B\odot\sF_{t'_B}$,
and since $\pi^B$ is universal, there exists a unique
natural transformation $\delta^B:t_B\Rightarrow t'_B$
such that $\fib_B\circ\Delta=\pi^B*\sF_{\delta^B}$. It remains to
check that the rule : $B\mapsto\delta^B$ yields a modification
$\delta:t\leadsto t'$, since in this case we get
$\Delta=\pi*\sF_\delta$, which will prove that
\eqref{eq_back-belt} is also full, thus concluding the
proof of the theorem. However, $\Delta$ is the datum of
a modification $\Delta_i:\pi_i\odot t\leadsto\pi_i\odot t'$
for every $i\in\Ob(I)$; the compatibility condition for $\Delta_i$
asserts that the following two compositions of oriented squares
coincide for every morphism $g:B'\to B$ of $\cB$ :
$$
\xymatrix@C+12pt{
X(B) \ar[r]^-{t_B} \ddouble \drtwocell\omit{_\ \ \delta^B} &
L(B) \ddouble \ar[r]^-{\pi^B_i}
\drtwocell\omit{_\ \ \ \ \ \ \one_{\pi^B_i}} & Fi(B) \ddouble &
X(B) \ar[r]^-{t_B} \ar[d]_{X(g)} \drtwocell\omit{_\ \ \ \tau^t_g} &
L(B) \ar[d]|{L(g)} \ar[r]^-{\pi^B_i}
\drtwocell\omit{_\ \ \ \ \tau^{\pi_i}_g} & Fi(B) \ar[d]^{Fi(g)} \\
X(B) \ar[r]|-{t'_B} \ar[d]_{X(g)} \drtwocell\omit{_\ \ \tau^{t'}_g} &
L(B) \ar[r]|-{\pi^B_i} \ar[d]|{L(g)}
\drtwocell\omit{_\ \ \ \ \tau^{\pi_i}_g} & Fi(B) \ar[d]^{Fi(g)} &
X(B') \ar[r]|-{t_{B'}} \ddouble \drtwocell\omit{_\ \ \delta^{B'}} &
L(B') \ar[r]|-{\pi^{B'}_i} \ddouble
\drtwocell\omit{_\ \ \ \ \ \ \one_{\pi^{B'}_i}} & Fi(B') \ddouble \\
X(B') \ar[r]_-{t'_{B'}} & L(B') \ar[r]_-{\pi^{B'}_i} & Fi(B') &
X(B') \ar[r]_-{t'_{B'}} & L(B') \ar[r]_-{\pi^{B'}_i} & Fi(B')
}$$
where $\tau^t_g$, $\tau^{t'}_g$ and $\tau^{\pi_i}_g$ are the
coherent constraints of respectively $t$, $t'$ and $\pi_i$.
We deduce that
$$
\pi^{B'}_i*(\tau^{t'}_g\odot(L(g)*\delta^B))=
\pi^{B'}_i*((\delta^{B'}*X(g))\odot\tau^t_g)
\qquad
\text{for every $i\in\Ob(I)$}
$$
and invoking again the universality of $\pi^{B'}$ we conclude
that $\tau^{t'}_g\odot(L(g)*\delta^B)=(\delta^{B'}*X(g))\odot\tau^t_g$,
which is the required compatibility condition for $\delta$.
\end{proof}

\subsection{Fibrations in groupoids}
\label{sec_fibrations-groupoids}
Recall that a {\em groupoid} is a category all whose morphisms
are invertible. To every category $\cB$, we may attach the
groupoid $\cB^\times$ such that $\Ob(\cB^\times)=\Ob(\cB)$,
and whose morphisms are the isomorphisms of $\cB$; the
composition law for morphisms in $\cB^\times$ is the same
as that of $\cB$, so $\cB^\times$ is a subcategory of $\cB$.
Then clearly every functor $G:\cA\to\cB$ restricts to
a functor $G^\times:\cA^\times\to\cB^\times$, and for every
universe $\sV$, the rules : $\cC\mapsto\cC^\times$ and
$(F:\cA\to\cB)\mapsto F^\times$ yield a right adjoint
$(-)^\times:\sV\tdu\bCat\to\sV\tdu\Gpd$ for the inclusion functor
$$
\sV\tdu\Gpd\to\sV\tdu\bCat
$$
of the full subcategory $\sV\tdu\Gpd$ of $\sV\tdu\bCat$
whose objects are the $\sV$-small groupoids.

\begin{remark}\label{rem_cat-groupoid}
(i)\ \
Notice that every isomorphism of functors $\alpha:G\isom H$
restricts to an isomorphism of functors
$\alpha^\times:G^\times\isom H^\times$.

(ii)\ \
Let $\cA$ and $\cB$ be any two categories, and $G:\cA\to\cB$
any functor. Then it is easily seen that $G$ is essentially
surjective if and only if the same holds for $G^\times$.
Moreover, if $G$ is $i$-faithful for some $i\in\{0,1,2\}$,
then the same holds for $G^\times$. The verification in case
$i=0$ shall be left to the reader. Next, let $f:GX\to GY$
be a morphism in $\cB^\times$; if $G$ is fully faithful,
there exists $g:X\to Y$ in $\cA$ such that $Gg=f$. But
$f$ is invertible, so by the same token there exists $h:Y\to X$
in $\cA$ such that $Gh=f^{-1}$; since $G(f\circ g)=G\one_Y$
and $G(g\circ f)=G\one_X$, we deduce $f\circ g=\one_Y$ and
$g\circ f=\one_X$, since $G$ is faithful. This shows that
$g$ is a morphism in $\cA^\times$, and proves the assertion
for $i=1$. For $i=2$, the functor $G$ is an equivalence,
so we know already that $G^\times$ is fully faithful; but
we have also already noticed that $G^\times$ is essentially
surjective, whence the assertion for $i=2$.

(iii)\ \
Since $(-)^\times$ is a right adjoint, it commutes with
all limits of $\sV\tdu\bCat$ (see proposition
\ref{prop_was-get-maddd}(iii)). Moreover, by inspecting
the construction of example \ref{ex_fil-colim-in-Cat},
it is easily seen that  $(-)^\times$ also commutes with
all filtered (small) colimits : details left to the reader.
\end{remark}

\sset\subsubsection{}\label{subsec_up-to-gpds-in-fibs}
We wish now to upgrade the associated groupoid construction
to the $2$-category of $\cC$-fibrations, for any given base
category $\cC$. Indeed, let $F:\cA\to\cC$ be any fibration;
we associate with $\cA$ the category
$$
\cA^\times
$$
with $\Ob(\cA^\times)=\cA$, and whose morphisms are the cartesian
morphisms of $\cA$. The composition law for morphisms in
$\cA^\times$ is the restriction of that of $\cA$, hence
$\cA^\times$ is a subcategory of $\cA$, and we denote by
$F^\times:\cA^\times\to\cC$ the restriction of $F$. Then it
is easily seen that $F^\times$ is also a fibration, and
all morphisms of $\cA^\times$ are cartesian for this fibration;
moreover, for every $X\in\Ob(\cC)$ the fibre category
$(F^\times)^{-1}X$ is a groupoid. The latter assertion can
be easily checked directly, and it follows also
straightforwardly from lemma \ref{lem_distinguished-cleavage}(i).
We call $F^\times$ the {\em fibration in groupoids associated
with\/} $F$. We say that $F$ is a {\em fibration in groupoids}
if $\cA=\cA^\times$. Clearly every cartesian functor $G:\cA\to\cB$
of $\cC$-fibrations restricts to a (cartesian) functor
$$
G^\times:\cA^\times\to\cB^\times.
$$
On the other hand, a given natural $\cC$-transformation
$\beta:G\Rightarrow H$ between $\cC$-cartesian functors
$G,H:\cA\to\cB$ induces a natural $\cC$-transformation
$\beta^\times:G^\times\Rightarrow H^\times$ if and only if
$\beta$ is an isomorphism of functors, {\em i.e.} if and
only if $\beta$ is a morphism of the groupoid
$\sCart_\cC(\cA,\cB)^\times$ associated with
$\sCart_\cC(\cA,\cB)$ as in \eqref{sec_fibrations-groupoids}.
Thus, we get a natural functor
\set\begin{equation}\label{eq_groupoids}
\sCart_\cC(\cA,\cB)^\times\to\sCart_\cC(\cA^\times,\cB^\times)
\qquad
G\mapsto G^\times
\qquad
(\beta:G\Rightarrow H)\mapsto\beta^\times
\end{equation}

\sset\subsubsection{}\label{subsec_Gpd}
For every universe $\sV$, let us also denote by
$$
\sV\tdu\Fib^\times(\cC)
\qquad\text{and}\qquad
\sV\tdu\Gpd(\cC)
$$
respectively : the subcategory of $\sV\tdu\Fib(\cC)$ whose objects
are all the $\cC$-fibrations with $\sV$-small fibres, and whose
$\Hom$-category is $\sCart_\cC(\cA,\cB)^\times$ for every pair
$(\cA,\cB)$ of such $\cC$-fibrations, and : the strong
sub-$2$-category of $\sV\tdu\Fib(\cC)$ whose objects are the
fibrations in groupoids. With this notation, clearly we get
as well a strict pseudo-functor
$$
(-)^\times_\cC:\sV\tdu\Fib^\times(\cC)\to\sV\tdu\Gpd(\cC)
\qquad
\cA\mapsto\cA^\times
$$
which is given on $\Hom$-categories by the foregoing system
of functors \eqref{eq_groupoids}.

\sset\subsubsection{}\label{subsec_times-and-upper-*}
Moreover, every functor $u:\cC\to\cC'$ induces strict
pseudo-functors
$$
\sV\tdu\Fib^\times(u)^*:
\sV\tdu\Fib^\times(\cC')\to\sV\tdu\Fib^\times(\cC)
\qquad
\sV\tdu\Gpd(u)^*:\sV\tdu\Gpd(\cC')\to\sV\tdu\Gpd(\cC)
$$
defined as the restrictions of
$\sV\tdu\Fib(u)^*:\sV\tdu\Fib(\cC')\to\sV\tdu\Fib(\cC)$, and by
simple inspection we get a commutative diagram of $2$-categories :
$$
\xymatrix@C+40pt{
\sV\tdu\Fib^\times(\cC') \ar[r]^-{\sV\tdu\Fib^\times(u)^*}
\ar[d]_{(-)^\times_{\cC'}} & \sV\tdu\Fib^\times(\cC) \ar[d]^{(-)^\times_\cC} \\
\sV\tdu\Gpd(\cC') \ar[r]^-{\sV\tdu\Gpd(u)^*} &
\sV\tdu\Gpd(\cC).
}$$

\sset\subsubsection{}\label{subsec_times-and-lower-*}
Let $u:\cC\to\cC'$ be a functor from a $\sV$-small category $\cC$
to a category $\cC'$ with $\sV$-small $\Hom$-sets. Recall that the
source functor $\ss_X:u\cC/X\to\cC$ is a fibration in groupoids
for every $X\in\Ob(\cC')$ (example \ref{ex_fibred-cats}(i));
taking into account remark \ref{rem_explicit-Fib_*}(ii) we
deduce for every $\cC$-fibration $\cA$ with $\sV$-small fibres
and every such $X$, natural equivalences of categories :
$$
(\sV\tdu\Fib(u)_*\cA(X))^\times\isom\sCart_\cC(u\cC/X,\cA)^\times
\isom\sCart_\cC(u\cC/X,\cA^\times)\isom(\sV\tdu\Fib(u)_*(\cA^\times))(X).
$$
There follows a pseudo-commutative diagram of $2$-categories :
$$
\xymatrix@C+40pt{
\sV\tdu\Fib^\times(\cC) \ar[r]^-{\sV\tdu\Fib^\times(u)_*}
\ar[d]_{(-)^\times_\cC} & \sV\tdu\Fib^\times(\cC') \ar[d]^{(-)^\times_{\cC'}} \\
\sV\tdu\Gpd(\cC) \ar[r]^-{\sV\tdu\Gpd(u)_*} &
\sV\tdu\Gpd(\cC')
}$$
where $\sV\tdu\Fib^\times(u)_*$ and $\sV\tdu\Gpd(u)_*$ are
the restrictions of $\sV\tdu\Fib(u)_*$. {\em Suppose
moreover that the category $X/u\cC$ is cofiltered for
every $X\in\Ob(\cC')$}; then, in view of remarks
\ref{rem_cat-groupoid}(iii) and
\ref{rem_alternative-inverse-image}(iii), for every such
$X$ we have as well a natural equivalence of categories :
$$
(\sV\tdu\Fib(u)_!\cA(X))^\times\isom\colim_{X\to uY}\cA(Y)^\times
\isom\sV\tdu\Fib(u)_!(\cA^\times)(X)
$$
where the colimit ranges over the small filtered category
$(X/u\cC)^o$. We get therefore a pseudo-commutative diagram
of $2$-categories :
\set\begin{equation}\label{eq_awdnews}
{\spreaddiagramcolumns{+40pt}\diagram
\sV\tdu\Fib^\times(\cC) \ar[r]^-{\sV\tdu\Fib^\times(u)_!}
\ar[d]_{(-)^\times_\cC} & \sV\tdu\Fib^\times(\cC') \ar[d]^{(-)^\times_{\cC'}} \\
\sV\tdu\Gpd(\cC) \ar[r]^-{\sV\tdu\Gpd(u)_!} &
\sV\tdu\Gpd(\cC)
\enddiagram}
\end{equation}
where $\sV\tdu\Fib^\times(u)_!$ and $\sV\tdu\Gpd(u)_!$
denote the restrictions of $\sV\tdu\Fib(u)_!$.

\begin{remark}\label{rem_groupoids}
Let $\cC$ be any category and $\sV$ any universe.

(i)\ \
It is easily seen that the pseudo-functor $(-)^\times_\cC$
is a strong right $2$-adjoint for the inclusion
pseudo-functor $\sV\tdu\Gpd(\cC)\to\sV\tdu\Fib^\times(\cC)$,
and the counit of the resulting $2$-adjunction assigns to
every fibration $\cA\to\cC$ the inclusion functor
$\cA^\times\to\cA$ : details left to the reader.

(ii)\ \
Notice that a fibration $\pi:\cA\to\cC$ is a fibration in
groupoids if and only if $\pi^{-1}X$ is a groupoid for every
$X\in\Ob(\cC)$. Indeed, the condition is obviously necessary.
Conversely, recall that every morphism $f$ of $\cA$ is the
composition of a cartesian morphism and a morphism
$g$ such that $\pi(g)=\one_X$ for some $X\in\Ob(\cC)$; but
if $\pi^{-1}X$ is a groupoid, then $g$ is an isomorphism
in $\cA$, in particular it is cartesian, so the same holds
for $f$, which shows that $\cA=\cA^\times$.

(iii)\ \
By the same token, to every pseudo-functor $F:\cC^o\to\bCat$
we attach the pseudo-functor
$$
F^\times:\cC^o\to\bCat
\qquad
X\mapsto(FX)^\times
$$
(whose coherence constraints are the same as those of $F$).
Then we easily see that :
$$
\cFib(F)^\times=\cFib(F^\times).
$$
\end{remark}

\sset\subsubsection{}\label{subsec_pi_0-on-groupoids}
With the notation of \eqref{subsec_Gpd}, notice that for
every presheaf $F\in\Ob(\cC^\wedge_\sV)$, the fibration
$\cFib_\cC(F)$ is a fibration in groupoid, so we get a
well defined strict pseudo-functor
$$
\cFib_\cC:\cC^\wedge_\sV\to\sV\tdu\Gpd(\cC)
$$
and by restriction, the pseudo-functor $\pi_0^\cC$
(see \eqref{subsec_connex-comp-fibration}) yields
a strong left $2$-adjoint
$$
\pi_0^\cC:\sV\tdu\Gpd(\cC)\to\cC^\wedge_\sV.
$$
Moreover, for every functor $u:\cC\to\cC'$, a simple
inspection yields an essentially commutative diagram
$$
\xymatrix@C+30pt{
\sV\tdu\Gpd(\cC') \ar[r]^-{\sV\tdu\Gpd(u)^*} \ar[d]_{\pi_0^{\cC'}} &
\sV\tdu\Gpd(\cC) \ar[d]^{\pi_0^\cC} \\
\cC'^\wedge_\sV \ar[r]^-{u^\wedge_\sV} & \cC^\wedge_\sV.
}$$
Furthermore, if $\cC$ is $\sV$-small and $\cC'$ has
$\sV$-small $\Hom$-sets, we also get the essentially
commutative diagram :
\set\begin{equation}\label{eq_lost-at-sea}
{\spreaddiagramcolumns{+30pt}\diagram
\sV\tdu\Gpd(\cC) \ar[r]^-{\sV\tdu\Gpd(u)_!} \ar[d]_{\pi_0^\cC} &
\sV\tdu\Gpd(\cC') \ar[d]^{\pi_0^{\cC'}} \\
\cC^\wedge_\sV \ar[r]^-{u_{\sV!}} & \cC'^\wedge_\sV.
\enddiagram}
\end{equation}
Indeed, for every $X\in\Ob(\cC')$ and every fibration $\cA$
with $\sV$-small fibres we have natural equivalences of
categories :
$$
(\pi_0^{\cC'}\sV\tdu\Gpd(u)_!\cA)(X)\isom
\pi_0((\Fib(u)_!\cA)(X))\isom\pi_0(\Pscolim{X\to uY}\cA(Y))
\isom\Pscolim{X\to uY}\pi_0(\cA(Y))
$$
since $\pi_0$ is a $2$-left adjoint pseudo-functor (proposition
\ref{prop_2-adjs-and-lims}(ii)). However, remark
\ref{rem_pseudo-limit}(vii) says that the $2$-colimit of
the (strict) pseudo-functor given by the rule
$(X\to uY)\mapsto\pi_0(\cA(Y))$ also represents the
colimit of the same functor, in the category of $\sV$-small
sets; on the other hand, the latter colimit is represented
by $u_{\sV!}(\pi^\cC_0\cA)(X)$, whence the contention.

\subsection{Sieves and descent theory}
\label{sec_sieves}
This section develops the basics of descent theory, in the
general framework of fibred categories.

\begin{definition}\label{def_sieve}
Let $\cB$ and $\cC$ be two categories, $F:\cB\to\cC$ a functor.
\begin{enumerate}
\item
A {\em sieve\/} of $\cC$ is a full subcategory $\cS$ of $\cC$ such
that the following holds. If $A\in\Ob(\cS)$, and $B\to A$ is any
morphism in $\cC$, then $B\in\Ob(\cS)$.
\item
If $S\subset\Ob(\cC)$ is any subset, there is a smallest sieve
$\cS_{\!S}$ of $\cC$ such that $S\subset\Ob(\cS_{\!S})$; we call
$\cS_{\!S}$ the {\em sieve generated by $S$}. If $S'\subset\Ob(\cC)$
is another subset and $\cS_{\!S'}\subset\cS_{\!S}$, we say that $S'$
is a {\em refinement} of $S$.
\item
If $\cS$ is a sieve of $\cC$, the
{\em inverse image of $\cS$ under $F$} is the full subcategory
$F^{-1}\!\cS$ of $\cB$ with
$\Ob(F^{-1}\!\cS)=\{B\in\Ob(\cB)~|~FB\in\Ob(\cS)\}$ (notice that
$F^{-1}\!\cS$ is a sieve).
\item
If $f:X\to Y$ is any morphism in $\cC$, and $\cS$ is any sieve of
$\cC\!/Y$, we shall write $\cS\times_Yf$ for the inverse image of
$\cS$ under the functor $f_*$ (notation of \eqref{eq_push-for}).
\end{enumerate}
\end{definition}

\begin{remark}\label{rem_sieves-and-sub}
Let $\cC$ be a category with small $\Hom$-sets and $X$
an object of\/ $\cC$.

(i)\ \
Every sieve $\cT$ of $\cC$ yields a subobject $F_\cT$
of the final object $\bone_\cC$ of $\cC^\wedge$, by
declaring that $F_\cT(Y)\neq\emptyset$ if and only if
$Y\in\Ob(\cT)$; hence $F_\cT(Y)$ is a set with one element
for every $Y\in\Ob(\cT)$, and is empty for
$Y\in\Ob(\cC)\setminus\Ob(\cT)$. It is
easily seen that the rule $\cT\mapsto F_\cT$ establishes
a bijection between the set of sieves of $\cC$ and the
set of subobjects of $\bone_\cC$. The inverse mapping
is given by the rule :
$$
(F\subset\bone_\cC)\mapsto(\cFib(F)\subset\cFib(\bone_\cC)=\cC).
$$

(ii)\ \
For every sieve $\cS$ of the category $\cC/X$ we
may then also consider the presheaf on $\cC$ 
$$
h_\cS:=\ss_{X!}(F_\cS)
$$
where $F_\cS$ is the presheaf on $\cC/X$ defined as in
(i), and $\ss_{X!}:(\cC/X)^\wedge\to\cC^\wedge$ is the left
adjoint to the functor $\ss^\wedge_X$ induced by the source
functor $\ss_X:\cC/X\to\cC$ (proposition
\ref{prop_in-the-same-vein}(vi.a)). By proposition
\ref{prop_in-the-same-vein}(vi.b,c), the presheaf $h_\cS$
is a subobject of $h_X$ (notation of \eqref{subsec_yoneda}),
and by inspecting \eqref{eq_can-be-used} we see that
$$
h_\cS(Y)=\{f\in\Hom_\cC(Y,X)~|~(Y,f)\in\Ob(\cS)\}
\qquad
\text{for every $Y\in\Ob(\cC)$}.
$$
For a given morphism $f:Y'\to Y$ in $\cC$, the map $h_\cS(f)$
is just the restriction of $\Hom_\cC(f,X)$. By the same token,
it also follows that the rule $\cS\mapsto h_\cS$ sets up a
natural bijection between the subobjects of $h_X$ in
$\cC^\wedge$ and the sieves of $\cC/X$. The inverse
mapping is given by the rule :
$$
(F\subset h_X)\mapsto(\cFib(F)\subset\cFib(h_X)=\cC/X)
$$
(see example \ref{ex_fibred-cats-II}(i)). Especially, the
restriction $\cS\to\cC$ of the source functor
$\ss_X:\cC/X\to\cC$ is a fibration, and example
\ref{ex_fibred-cats-II}(i) generalizes to a natural
isomorphism of $\cC$-fibrations :
$$
\cFib(h_\cS)\isom\cS.
$$
Combining with lemma \ref{lem_lable}, we also deduce
a natural isomorphism in $\cC^\wedge$ :
\set\begin{equation}\label{eq_colim-sieve}
\colim_\cS h_\cC\circ\ss_\cS\isom h_\cS
\end{equation}
where $\ss_\cS:\cS\to\cC$ is the restriction of the source
functor $\ss_X$ of \eqref{subsec_slice-cat}.

(iii)\ \
Let $S:=\{X_i\to X~|~i\in I\}$ be any family of morphisms
of\/ $\cC$. Then $S$ generates a given sieve $\cS$ of
$\cC/X$ if and only if :
$$
h_\cS=\bigcup_{i\in I}\Img(h_{X_i}\to h_X).
$$
(Notice that the above union is well defined even in case $I$
is not small.)

(iv)\ \
Moreover, if $S$ as in (iii) generates $\cS$ and $f:Y\to X$ is
any morphism such that the fibre product $Y_i:=Y\times_XX_i$
is representable in $\cC$ for every $i\in I$, then
$\cS\times_Xf$ is the sieve generated by the family of induced
projections $\{Y_i\to Y~|~i\in I\}\subset\Ob(\cC\!/\!Y)$.

(v)\ \
Furthermore, let $f:Y\to X$ be any morphism of $\cC$; it is
easily seen that the correspondence of (ii) induces a natural
identification of subobjects of $h_Y$:
$$
h_{\cS\times_Xf}=h_\cS\times_{h_X} h_Y
\qquad
\text{for every sieve $\cS$ of $\cC/X$}.
$$
Since the Yoneda embedding is fully faithful, we may
sometimes abuse notation, to identify $X$ with the
corresponding representable presheaf $h_X$; then we may
write $h_{\cS\times_Xf}=h_\cS\times_XY$. In view of (ii)
we deduce an isomorphism of fibred $\cC$-categories
$\cFib(h_\cS\times_XY)\isom\cS\times_Xf$.
On the other hand, let $\pi:h_\cS\times_XY\to h_\cS$ be
the natural projection; a direct inspection shows that the
foregoing isomorphisms fit into the commutative diagram :
$$
\xymatrix{ \cFib(h_\cS\times_XY) \ar[r]^-\sim
\ar[d]_{\cFib(\pi)} & \cS\times_Xf \ar[d]^{f_{*|\cS}} \\
\cFib(h_\cS) \ar[r]^-\sim & \cS
}$$
where the right vertical arrow is the restriction of the
functor $f_*:\cC/Y\to\cC/X$ (details left to the reader).
(Especially, $f_{*|\cS}$ is a $\cC$-cartesian functor, but
this assertion is trivial, since every morphism of $\cC/X$
is $\cC$-cartesian).
\end{remark}

\sset\subsubsection{}\label{subsec_step-one}
Let $\cC$ be a small category, and $\cS$ the sieve of $\cC$
generated by a subset $S\subset\Ob(\cC)$. Say that
$S=\{S_i~|~i\in I\}$ for a small set $I$; for every $i\in I$
there is a faithful embedding $\eps_i:\cC\!/\!S_i\to\cS$, and
for every pair $(i,j)\in I\times I$, we define
$$
\cC\!/\!S_{ij}:=\cC\!/\!S_i\times_{(\eps_i,\eps_j)}\cC\!/\!S_j
$$
(notation of example \ref{ex_cat-cats}(i)). Hence,
the objects of $\cC\!/\!S_{ij}$ are all the triples
$(X,g_i,g_j)$, where $X\in\Ob(\cC)$ and $g_l\in\Hom_\cC(X,S_l)$
for $l=i,j$. The natural projections :
$$
\pi^1_{ij*}:\cC\!/\!S_{ij}\to\cC\!/\!S_i
\qquad
\pi^0_{ij*}:\cC\!/\!S_{ij}\to\cC\!/\!S_j
$$
are faithful embeddings. We deduce a natural diagram of categories :
$$
\coprod_{(i,j)\in I\times I}\!\!\!\!\!\!\! \xymatrix{\cC\!/\!S_{ij}
\ar@<.5ex>[r]^-{\partial_0} \ar@<-.5ex>[r]_-{\partial_1} &}
\!\!\coprod_{i\in I}\cC\!/\!S_i \xrightarrow{\eps} \cS
$$
where :
$$
\partial_0:=\!\!\!\!\!\coprod_{(i,j)\in I\times I}\!\!\!\pi^0_{ij*}
\qquad
\partial_1:=\!\!\!\!\!\coprod_{(i,j)\in I\times I}\!\!\!\pi^1_{ij*}
\qquad \eps:=\coprod_{i\in I}\eps_i.
$$

\begin{remark}\label{rem_abuse-already}
Notice that, under the current assumptions,
the product $S_{ij}:=S_i\times S_j$ is not necessarily
representable in $\cC$. In case it is, we may consider
another category, also denoted $\cC\!/\!S_{ij}$, namely
the category of $S_{ij}$-objects of $\cC$ (as in
\eqref{subsec_slice-cat}). The latter is naturally
isomorphic to the category with the same name introduced
in \eqref{subsec_step-one}. Moreover, under this natural
isomorphism, the projections $\pi^0_{ij*}$ and $\pi^1_{ij*}$
are identified with the functors induced by the natural
morphisms $\pi^0_{ij}:S_{ij}\to S_j$ and respectively
$\pi^1_{ij}:S_{ij}\to S_i$. Hence, in this case, the
notation of \eqref{subsec_step-one} is compatible with
\eqref{eq_push-for}.
\end{remark}

\begin{lemma}\label{lem_coeq-sieve}
With the notation of \eqref{subsec_step-one}, the functor $\eps$
induces an isomorphism between $\cS$ and the coequalizer (in the
category $\bCat$) of the pair of functors $(\partial_0,\partial_1)$.
\end{lemma}
\begin{proof} Let $\cA$ be any other object of $\bCat$, and
$F:\coprod_{i\in I}\cC\!/\!S_i\to\cA$ a functor such that
$F\circ\partial_0=F\circ\partial_1$. We have to show that $F$
factors uniquely through $\eps$. To this aim, we construct
explicitly a functor $G:\cS\to\cA$ such that $G\circ\eps=F$. First
of all, by the universal property of the coproduct, $F$ is the same
as a family of functors $(F_i:\cC\!/\!S_i\to\cA~|~i\in I)$, and the
assumption on $F$ amounts to the system of identities :
\set\begin{equation}\label{eq_amounts-to}
F_i\circ\pi^1_{ij*}=F_j\circ\pi^0_{ij*} \qquad\text{for every
$i,j\in I$}.
\end{equation}
Hence, let $X\in\Ob(\cS)$; by assumption there exist $i\in I$ and a
morphism $f:X\to S_i$ in $\cC$, so we may set $GX:=F_if$. In case
$g:X\to S_j$ is another morphism in $\cC$, we deduce an object
$h:=(X,f,g)\in\Ob(\cC\!/\!S_{ij})$, so $f=\pi^1_{ij*}h$ and
$g=\pi^0_{ij*}h$; then \eqref{eq_amounts-to} shows that $F_if=F_jg$,
{\em i.e.} $GX$ is well-defined.

Next, let $\phi:X\to Y$ be any morphism in $\cS$; choose $i\in I$
and a morphism $f_Y:Y\to S_i$, and set $f_X:=f_Y\circ\phi$. We let
$G\phi:=F_i(\phi:f_X\to f_Y)$. Arguing as in the foregoing, one
verifies easily that $G\phi$ is independent of all the choices, and
then it follows easily that $G(\psi\circ\phi)=G\psi\circ G\phi$ for
every other morphism $\psi:Y\to Z$ in $\cS$. It is also clear that
$G\one_X=\one_{GX}$, whence the contention.
\end{proof}

\sset\subsubsection{}\label{subsec_towards-top}
In the situation of \eqref{subsec_step-one}, suppose that the
set of generators $S$ is the whole of $\Ob(\cS)$; in this case,
the augmentation $\eps$ can also be used to produce the following
$2$-categorical presentation of $\cS$, which upgrades the
isomorphism \eqref{eq_colim-sieve}. Consider the strict
pseudo-functor
$$
G_\cS:\cS\to\Fib(\cC)
\quad : \quad
Y\mapsto\cC\!/Y
\qquad
(Z\xrightarrow{\ f\ }Y)\mapsto(\cC/Z\xrightarrow{\ f_*\ }\cC/Y).
$$
We have a natural strict pseudo-cocone
$$
\hat\eps:G_\cS\Rightarrow\sF_\cS
$$
where $\sF_\cS$ is the constant pseudo-functor $\cS\to\Fib(\cC)$ with
value $\cS$, and $\hat\eps_X:\cC/X\to\cS$ is the faithful embedding
as in \eqref{subsec_step-one}, for every $X\in\Ob(\cS)$. We may
then state :

\begin{lemma}\label{lem_new-pseudo-col}
The pseudo-cocone $\hat\eps$ induces an equivalence of fibrations
over $\cC$ :
$$
\Pscolim{\cS}G_\cS\isom\cS.
$$
\end{lemma}
\begin{proof} By theorem \ref{th_fib-limits-are-fibrewise}(ii),
the $2$-colimits are computed fibrewise in the $2$-category
$\Fib(\cC)$, so we are reduced to checking that $\widehat\eps$
induces an equivalence of categories
$$
\Pscolim{\cS}G_{\cS,X}\isom\cS_X
\qquad
\text{for every $X\in\Ob(\cC)$}
$$
where $\cS_X$ is the fibre category over $X$ of the fibration
$\cS\to\cC$, and $G_{\cS,X}:\cS\to\bCat$ is the strict
pseudo-functor that assigns to every $Y\in\Ob(\cS)$ the
fibre category $(\cC/Y)_X$ of the fibration $\cC/Y\to\cC$
(for every morphism $f:Z\to Y$ in $\cS$, the corresponding
functor $G_{\cS,X}(f)$ is the restriction of $G_\cS(f)=f_*$).
However, if $X\notin\Ob(\cS)$, the category $\cS_X$ is empty,
and $G_{\cS,X}$ is the constant pseudo-functor with value
equal to the empty category, so the assertion is clear in
this case. Suppose then that $X\in\Ob(\cS)$, in which case
$\cS_X$ is the category whose unique object is $X$ and whose
unique morphism is $\one_X$. Also, for every $Y\in\Ob(\cS)$,
the category $(\cC/Y)_X$ is discrete, with set of objects
given by $\Hom_\cC(X,Y)$. According to example
\ref{ex_filter-2-colim-in-Cat}(i), the strong $2$-colimit
of $G_{\cS,X}$ is represented by the category
$\cFib(G_{\cS,X})[\Sigma^{-1}]$, where $\Sigma$ is the set of
cartesian morphisms of the fibration $\cFib(G_{\cS,X})\to\cS^o$.
Now, it is easily seen that $\cFib(G_{\cS,X})=(X/\cS)^o=\cS^o/X^o$,
with structure functor given by the usual source functor
$\cS^o/X^o\to\cS^o$. Hence, $\Sigma$ is the set of all morphisms
of $\cC^o/X^o$, and since $\cC^o/X^o$ has the final object
$\one_{X^o}$, the assertion follows easily from example
\ref{ex_invert-all-with-final}.
\end{proof}

\begin{definition}\label{def_descent-fibred}
Let $\phi:\cA\to\cB$ be a fibration, $B\in\Ob(\cB)$ an object,
$\cS\subset\cB/B$ a sieve, $i\in\{0,1,2\}$, and denote by
$\iota_\cS:\cS\to\cB/B$ the fully faithful embedding.
\begin{enumerate}
\item
We say that $\cS$ is a {\em sieve of $\phi$-$i$-descent},
if the restriction functor :
$$
\sCart_\cB(\iota_\cS,\cA):\cA(B)\to\sCart_\cB(\cS,\cA)
$$
is $i$-faithful (see remark \ref{rem_i-faithful}).
\item
We say that $\cS$ is a {\em sieve of universal $\phi$-$i$-descent}
if the sieve $\cS\times_Bf$ is of $\phi$-$i$-descent for every
morphism $f:B'\to B$ of $\cB$ (notation of definition
\ref{def_sieve}(iv)).
\item
Let $f:B'\to B$ be a morphism in $\cB$. We say that $f$ is a
{\em morphism of $\phi$-$i$-descent} (resp. a {\em morphism of
universal $\phi$-$i$-descent}), if the sieve generated by $\{f\}$
is of $\phi$-$i$-descent (resp. of universal $\phi$-$i$-descent).
\end{enumerate}
\end{definition}

\begin{remark}\label{rem_fin-prod-reps}
In the situation of definition \ref{def_descent-fibred},
suppose that $\cS$ is the sieve generated by a set of
objects $\{S_i\to B~|~i\in I\}\subset\Ob(\cB/B)$.
There follows a natural diagram of categories (notation
of \eqref{subsec_step-one}) :
$$
\sCart_\cB(\cS,\cA)\xrightarrow{\eps^*}\prod_{i\in I}
\xymatrix{\cA(S_i)
\ar@<.5ex>[r]^-{\partial_0^*} \ar@<-.5ex>[r]_-{\partial_1^*} &}
\!\!\!\!\!\prod_{(i,j)\in I\times I}\!\!\!\sCart_\cB(\cB/S_{ij},\cA)
$$
where $\eps^*:=\sCart_\cB(\eps,\cA)$ and
$\partial_i^*:=\sCart_\cB(\partial_i,\cA)$, for $i=0,1$.
With this notation, lemma \ref{lem_coeq-sieve} easily
implies that $\eps^*$ induces an isomorphism between
$\sCart_\cB(\cS,\cA)$ and the equalizer (in the category
$\bCat$) of the pair of functors $(\partial^*_0,\partial^*_1)$.
\end{remark}

\sset\subsubsection{}\label{subsec_triple-prods}
We would like to exploit the presentation of
$\sCart_\cB(\cS,\cA)$ in remark \ref{rem_fin-prod-reps},
in order to translate definition \ref{def_descent-fibred}
in terms of the fibre categories $\cA_{S_i}$ and
$\cA_{S_{ij}}$. The problem is that such a translation
must be carried out via a pseudo-natural equivalence
(namely $\sev$), and such equivalences do not respect
a presentation as above in terms of equalizers in the
category $\bCat$. What we need is to upgrade our
presentation of $\cS$ to a new one, which is preserved
by pseudo-natural transformations. This is achieved as
follows. Resume the general situation of \eqref{subsec_step-one}.
For every $i,j,k\in I$, set $\cC/S_{ijk}:=\cC/S_{ij}\times_\cC\cC/S_k$.
We have a natural diagram of categories :
\set\begin{equation}\label{eq_wecanview}
\coprod_{(i,j,k)\in I^3}\!\!\!\!\!\!
\xymatrix@C+10pt{\cC/S_{ijk} \ar@<1ex>[r]^-{\partial'_0}
\ar@<-1ex>[r]_-{\partial'_2} \ar[r]|-{\:\partial'_1} &}
\!\!\!\!\coprod_{(i,j)\in I^2}\!\!\!\!
\xymatrix{\cC/S_{ij}
\ar@<.5ex>[r]^-{\partial_0} \ar@<-.5ex>[r]_-{\partial_1} &}
\!\!\coprod_{i\in I}\cC/S_i \xrightarrow{\eps} \cS
\end{equation}
where $\partial'_0$ is the coproduct of the natural projections
$\pi^0_{ijk*}:\cC/S_{ijk}\to\cC/S_{jk}$ for every $i,j,k\in I$,
and likewise $\partial'_1$ (resp. $\partial'_2$) is the coproduct
of the projections $\pi^1_{ijk*}:\cC/S_{ijk}\to\cC/S_{ik}$ (resp.
$\pi^2_{ijk*}:\cC/S_{ijk}\to\cC/S_{ij}$). We can view
\eqref{eq_wecanview} as an
{\em augmented $2$-truncated semi-simplicial object} in
$\Fib(\cC)$, {\em i.e.} a functor :
$$
F_\cS:\Sigma_2^{+o}\to\Fib(\cC)
$$
from the opposite of the category $\Sigma^+_2$ whose
objects are the ordered sets $[-1]$, $[0]$, $[1]$
and $[2]$, and whose morphisms are the non-decreasing
injective maps (this is a subcategory of the category
$\Delta^{\!\wedge}_2$ of definition \ref{def_simplicial-cats}(iii)).

\begin{remark}\label{rem_abuse-not}
Suppose that finite products are representable in $\cB$,
and for every $i,j,k\in I$, set $S_{ij}:=S_i\times S_j$,
and $S_{ijk}:=S_{ij}\times S_k$. Just as in remark
\ref{rem_abuse-already}, the category $\cC/S_{ijk}$ of
$S_{ijk}$-objects of $\cC$ is naturally isomorphic to
the category with the same name introduced in
\eqref{subsec_triple-prods}, and under this isomorphism,
the functors $\pi^0_{ijk*}$ are identified with the functors
arising from the natural projections
$\pi^0_{ijk}:S_{ijk}\to S_{jk}$ (and likewise for $\pi^1_{ijk*}$
and $\pi^2_{ijk*}$).
\end{remark}

With this notation, denote by $\Sigma_2$ the full subcategory
of $\Sigma^+_2$ whose objects are the non-empty sets; we
have the following $2$-category analogue of lemma
\ref{lem_coeq-sieve} :

\begin{proposition}\label{prop_wecanview}
The augmentation $\eps$ induces an equivalence of categories :
$$
\Pscolim{\Sigma_2^o}F_\cS\isom\cS.
$$
\end{proposition}
\begin{proof} By theorem \ref{th_fib-limits-are-fibrewise}(ii),
the $2$-colimits are computed fibrewise in the $2$-category
$\Fib(\cC)$, so we are reduced to checking that $\eps$
induces an equivalence of categories
$$
\Pscolim{\Sigma^o_2}F_{\cS,X}\isom\cS_X
\qquad
\text{for every $X\in\Ob(\cC)$}
$$
where $\cS_X$ is the fibre category over $X$ of the fibration
$\cS\to\cC$, and $F_{\cS,X}:=\fib_X\circ F_\cS$, where
$\fib_X:\Fib(\cC)\to\bCat$ is defined as in
\eqref{subsec_lims-in-Fib}. Now, if $X\notin\Ob(\cS)$, the
category $\cS_X$ is empty, and $\fib_X\circ F_\cS$ is the
constant pseudo-functor with value equal to the empty category,
so the assertion trivially holds in this case. Suppose then
that $X\in\Ob(\cS)$, in which case $\cS_X$ is the category whose
unique object is $X$ and whose unique morphism is $\one_X$. Set
$$
T:=\coprod_{i\in I}\Hom_\cC(X,S_i)
$$
and notice that $F_{\cS,J}[0],F_{\cS,J}[1]$ and $F_{\cS,J}[2]$
are the discrete categories whose sets of objects are
respectively $T$, $T\times T$ and $T\times T\times T$.
Morever, the strict pseudo-functor $F_{\cS,X}$ corresponds
to the $2$-truncated semi-simplicial diagram of sets :
\set\begin{equation}\label{eq_this-is-T}
\xymatrix@C+10pt{T\times T\times T \ar@<1ex>[r]^-{\partial'_0}
\ar@<-1ex>[r]_-{\partial'_2} \ar[r]|-{\:\partial'_1} & T\times T
\ar@<.5ex>[r]^-{\partial_0} \ar@<-.5ex>[r]_-{\partial_1} & T}
\end{equation}
whose maps $\partial_i$ and $\partial'_i$ are the natural
projections. To conclude, it therefore suffices to show :

\begin{claim} For every set $T\neq\emptyset$, the $2$-colimit
in the $2$-category $\bCat$ of the system of discrete categories
\eqref{eq_this-is-T} is represented by the discrete category
with one object.
\end{claim}
\begin{pfclaim}[] The diagram \eqref{eq_this-is-T} can
be regarded as a presheaf $T_\bullet$ on the category
$\Sigma_2$, and according to example
\ref{ex_filter-2-colim-in-Cat}(i), the strong $2$-colimit
of $T_\bullet$ is represented by the category
$\cT:=\cFib(T_\bullet)[\Lambda^{-1}]$, where $\Lambda$ is
the set of cartesian morphisms of the fibration
$\cFib(T_\bullet)\to\Sigma_2$; but every morphism of
$\cFib(T_\bullet)$ is cartesian (see
\eqref{subsec_fibred-cats-II}). Denote by $\cT_0$ the category
with $\Ob(\cT_0)=T$, and such that $\Hom_{\cT_0}(t,t')$ contains
a unique element $\tau_{t,t'}$, for every $t,t'\in T$. We define
a functor $q:\cT_0\to\cT$ as follows. Recall that the objects
of $\cFib(T_\bullet)$ are the pairs $(j,\underline t)$, where
$j\in\{0,1,2\}$ and $\underline t\in T^{j+1}$. For every
$\underline t:=(t_0,t_1)\in T\times T$ we have morphisms
$(0,t_1)\xleftarrow{\delta_{0,\underline t}}(1,\underline t)
\xrightarrow{\delta_{1,\underline t}}(0,t_0)$ in $\cFib(T_\bullet)$,
and we set
$$
q(t_0):=(0,t_0)
\quad\text{and}\quad
q(\tau_{t_0,t_1}):=
[\delta_{0,(t_0,t_1)}]\circ[\delta_{1,(t_0,t_1)}]^{-1}
\qquad
\text{for every $t_0,t_1\in T$}
$$
where we denote by $[f]$ the image in $\cT$ of every morphism
$f$ of $\cFib(T_\bullet)$. Indeed, for every
$\underline t:=(t_0,t_1,t_2)\in T^3$ and $j=0,1,2$, we have as
well morphisms $\delta'_{j,\underline t}:(2,\underline t)\to
(1,\partial'_j(\underline t))$ in $\cFib(T_\bullet)$, and we may
compute :
$$
\begin{aligned}
q(\tau_{t_1,t_2})\circ q(\tau_{t_0,t_1})&\,=
[\delta_{0,\partial'_0\underline t}]\circ[\delta_{1,\partial'_0\underline t}]^{-1}
\circ[\delta_{0,\partial'_2\underline t}]\circ[\delta_{1,\partial'_2\underline t}]^{-1} \\
&\,=[\delta_{0,\partial'_0\underline t}]\circ[\delta'_{0,\underline t}]\circ
[\delta'_{0,\underline t}]^{-1}\circ[\delta_{1,\partial'_0\underline t}]^{-1}
\circ[\delta_{0,\partial'_2\underline t}]\circ[\delta'_{2,\underline t}]\circ
[\delta'_{2,\underline t}]^{-1}\circ[\delta_{1,\partial'_2\underline t}]^{-1} \\
&\,=[\delta_{0,\partial'_0\underline t}]\circ[\delta'_{0,\underline t}]\circ
[\delta'_{0,\underline t}]^{-1}\circ[\delta_{1,\partial'_0\underline t}]^{-1}
\circ[\delta_{1,\partial'_0\underline t}]\circ[\delta'_{0,\underline t}]\circ
[\delta'_{1,\underline t}]^{-1}\circ[\delta_{1,\partial'_1\underline t}]^{-1} \\
&\,=[\delta_{0,\partial'_0\underline t}]\circ[\delta'_{0,\underline t}]\circ
[\delta'_{1,\underline t}]^{-1}\circ[\delta_{1,\partial'_1\underline t}]^{-1} \\
&\,=[\delta_{0,\partial'_1\underline t}]\circ[\delta'_{1,\underline t}]\circ
[\delta'_{1,\underline t}]^{-1}\circ[\delta_{1,\partial'_1\underline t}]^{-1} \\
&\,=q(\tau_{t_0,t_2})
\end{aligned}
$$
as required. Next, let also $q':\cT\to\cT_0$ be the unique functor
given by the rules :
$$
(0,t_0)\mapsto t_0
\qquad
(1,(t_0,t_1))\mapsto t_0
\qquad
(2,(t_0,t_1,t_2))\mapsto t_0
\qquad
\text{for every $t_0,t_1,t_2\in T$}.
$$
Obviously $q'\circ q=\one_{\cT_0}$. Let $L:\cFib(T_\bullet)\to\cT$
be the localization; to conclude, it suffices to exhibit an
isomorphism of functors $\one_\cT\isom q\circ q'$, and corollary
\ref{cor_of-localization} further reduces to exhibiting an
isomorphism $\omega:L\isom q\circ q'\circ L$. We define $\omega$
by the rules :
$$
(0,t_0)\mapsto\one_{(0,t_0)}
\qquad
(1,(t_0,t_1))\mapsto[\delta_{1,(t_0,t_1)}]
\qquad
(2,(t_0,t_1,t_2))\mapsto[\delta_{1,(t_0,t_2)}\circ\delta'_{2,(t_0,t_1,t_2)}]
$$
for every $t_0,t_1,t_2\in T$. The naturality of $\omega$ follows
by a straightforward verification.
\end{pfclaim}
\end{proof}

\sset\subsubsection{}\label{subsec_desc-cats}
Resume the situation of remark \ref{rem_fin-prod-reps}, and
notice that all the categories appearing in \eqref{eq_wecanview}
are fibred over $\cB$ : indeed, every morphism in each of
these categories is cartesian, hence all the functors
appearing in \eqref{eq_wecanview} are cartesian.
Let us consider now the functor :
$$
\sCart_\cB(-,\cA):(\bCat/\cB)^o\to\bCat
$$
that assigns to every $\cB$-category $\cC$ the category
$\sCart_\cB(\cC,\cA)$. With the notation of remark
\ref{rem_fin-prod-reps}, we deduce a functor :
$$
\sCart_\cB(F_\cS,\cA):\Sigma^+_2\to\bCat
$$
and in light of the foregoing observations, proposition
\ref{prop_wecanview} easily implies that $\sCart_\cB(\eps,\cA)$
induces an equivalence of categories :
\set\begin{equation}\label{eq_onthe-basis}
\Pslim{\Sigma_2}\sCart_\cB(F_\cS,\cA)\isom\sCart_\cB(\cS,\cA).
\end{equation}
Next, suppose that the fibre products $S_{ij}:=S_i\times S_j$
and $S_{ijk}:=S_{ij}\times S_k$ are representable in $\cB$ (see
remark \ref{rem_abuse-not}); in this case, we may compose with
the pseudo-equivalence $\sev^\cA_{|\bullet}$ of remark
\ref{rem_mathsfev}(i) : combining with lemma
\ref{lem_pseudo-trivial} we finally obtain an equivalence
between the category $\sCart_\cB(\cS,\cA)$,
and the $2$-limit of the pseudo-functor
$\sd:=\sev\circ\sCart_\cB(F_\cS,\cA):\Sigma_2\to\bCat$ :
$$
\prod_{i\in I}\xymatrix{\cA_{S_i}
\ar@<.5ex>[r]^-{\partial^0} \ar@<-.5ex>[r]_-{\partial^1} &}
\!\!\!\!\prod_{(i,j)\in I^2}\!\!\!
\xymatrix{\cA_{S_{ij}} \ar@<1ex>[r]^-{\partial^0}
\ar@<-1ex>[r]_-{\partial^2} \ar[r]|-{\:\partial^1} &}
\!\!\!\!\prod_{(i,j,k)\in I^3}\!\!\cA_{S_{ijk}}.
$$
Of course, the coface operators $\partial^s$ on
$\prod_{i\in I}\cA_{S_i}$ decompose as products
of pull-back functors:
$$
\pi^{0*}_{ij}:\cA_{S_j}\to\cA_{S_{ij}}
\qquad
\pi^{1*}_{ij}:\cA_{S_i}\to\cA_{S_{ij}}
$$
attached -- via the chosen cleavage $\sc$ of $\phi$ --
to the projections $\pi^0_{ij}:S_{ij}\to S_j$ and
$\pi^1_{ij}:S_{ij}\to S_i$ (and likewise for the components
$\pi^{t*}_{ijk}$ of the other coface operators).

\sset\subsubsection{}\label{subsec_handier}
By inspecting the proof of theorem \ref{th_bCat-cpt},
we may give the following explicit description of this
$2$-limit. Namely, it is the category whose objects
are the data
$$
\underline X:=
(X_i,X_{ij},X_{ijk},
\xi^u_i,\xi^s_{ij},\xi^t_{ijk}~|~i,j,k\in I; s\in\{0,1\}; u,t\in\{0,1,2\})
$$
where :
$$
X_i\in\Ob(\cA_{S_i})\quad
X_{ij}\in\Ob(\cA_{S_{ij}})\quad
X_{ijk}\in\Ob(\cA_{S_{ijk}})\qquad\text{for every $i,j,k\in I$}
$$
and for every $i,j,k\in I$ :
$$
\begin{aligned}
\xi^0_i & :(\pi^1_{ik}\pi^1_{ijk})^*X_i\isom X_{ijk} \qquad &
\xi^0_{ij} & :\pi^{0*}_{ij}X_j\isom X_{ij} \qquad &
\xi^0_{ijk} & :\pi^{0*}_{ijk}X_{jk}\isom X_{ijk} \\
\xi^1_j & :(\pi^1_{jk}\pi^0_{ijk})^*X_j\isom X_{ijk} \qquad &
\xi^1_{ij} & :\pi^{1*}_{ij}X_i\isom X_{ij} \qquad &
\xi^1_{ijk} & :\pi^{1*}_{ijk}X_{ik}\isom X_{ijk} \\
\xi^2_k & :(\pi^0_{jk}\pi^0_{ijk})^*X_k\isom X_{ijk} \qquad &
& & \xi^2_{ijk} & :\pi^{2*}_{ijk}X_{ij}\isom X_{ijk}
\end{aligned}
$$
are isomorphisms related by the cosimplicial identities :
$$
\begin{aligned}
\xi^0_{ijk}\circ\pi^{0*}_{ijk}\xi^0_{jk} &
=\xi^2_k\circ\gamma^{00}_X \qquad &
\xi^1_{ijk}\circ\pi^{1*}_{ijk}\xi^0_{ik} &
=\xi^2_k\circ\gamma^{01}_X \\
\xi^2_{ijk}\circ\pi^{2*}_{ijk}\xi^0_{ij} &
=\xi^1_j\circ\gamma^{02}_X \qquad &
\xi^0_{ijk}\circ\pi^{0*}_{ijk}\xi^1_{jk} &
=\xi^1_j\circ\gamma^{10}_X \\
\xi^2_{ijk}\circ\pi^{2*}_{ijk}\xi^1_{ij} &
=\xi^0_i\circ\gamma^{12}_X \qquad &
\xi^1_{ijk}\circ\pi^{1*}_{ijk}\xi^1_{ik} &
=\xi^0_i\circ\gamma^{11}_X
\end{aligned}
$$
where $\gamma^{st}:=\gamma_{\sd(\partial^s),\sd(\partial^t)}:
\sd(\partial^s)\circ\sd(\partial^t)\Rightarrow\sd(\partial^s\circ\partial^t)$
denotes the coherence constraint of the cleavage $\sc$, for any
pair of arrows $(\partial^s,\partial^t)$ in the category $\Sigma_2$.
The morphisms $\underline X\to\underline Y$ in this category
are the systems of morphisms :
$$
(X_i\to Y_i,X_{ij}\to Y_{ij},X_{ijk}\to Y_{ijk}~|~i,j,k\in I)
$$
that are compatible in the obvious way with the various
isomorphisms. However, one may argue as in the proof of
proposition \ref{prop_wecanview}, to replace this category
by an equivalent one which admits a handier description :
given a datum $\underline X$, one can make up an isomorphic
datum $\underline
X^*:=(X_i,X^*_{ij},X^*_{ijk},\eta_i,\eta_{ij},\eta_{ijk})$,
by the rule :
$$
X_{ij}^*:=\pi^{1*}_{ij}X_i \qquad
X_{ijk}^*:=(\pi^1_{ik}\circ\pi^1_{ijk})^*X_i
$$
$$
\begin{aligned}
\eta^0_i & :=\one \qquad & \eta^1_{ij} & :=\one \qquad &
\eta^0_{ijk} & :=\eta^1_j\circ\gamma^{10}_X \\
\eta^1_j & :=(\xi^0_i)^{-1}\circ\xi^1_j \qquad & \eta^0_{ij} &
:=(\xi^1_{ij})^{-1}\circ\xi^0_{ij} \qquad &
\eta^1_{ijk} & :=\gamma^{11}_X \\
\eta^2_k & :=(\xi^0_i)^{-1}\circ\xi^2_k \qquad & & \qquad  &
\eta^2_{ijk} & :=\gamma^{12}_X.
\end{aligned}
$$
The cosimplicial identities for this new object are subsumed into a
single cocycle identity for $\omega_{ij}:=\eta^0_{ij}$. Summing up,
we arrive at the following description of our $2$-limit :
\begin{itemize}
\item
The objects are all the systems
$\underline X:=(X_i,\omega^X_{ij}~|~i,j\in I)$
where $X_i\in\Ob(\cA_{S_i})$ for every $i\in I$, and
$$
\omega^X_{ij}:\pi^{0*}_{ij}X_j\isom\pi^{1*}_{ij}X_i
$$
is an isomorphism in $\cA_{S_{ij}}$, for every $i,j\in I$,
fulfilling the cocycle identity :
$$
p^{2*}_{ijk}\omega^X_{ij}\circ
p^{0*}_{ijk}\omega^X_{jk}=
p^{1*}_{ijk}\omega^X_{ik}
\qquad\text{for every $i,j,k\in I$}
$$
where, for every $i,j,k\in I$, and $t=0,1,2$ we have set :
$$
p^{t*}_{ijk}\omega^X_{\bullet\bullet}:=
\gamma^{1t}_X\circ\pi^{t*}_{ijk}\omega^X_{\bullet\bullet}\circ
(\gamma^{0t}_X)^{-1}.
$$
\item
The morphisms $\underline X\to\underline Y$ are the
systems of morphisms $(f_i:X_i\to Y_i~|~i\in I)$ with :
\set\begin{equation}\label{eq_morphius}
\omega^Y_{ij}\circ\pi^{0*}_{ij}f_j=\pi^{1*}_{ij}f_i\circ\omega^X_{ij}
\qquad\text{for every $i,j\in I$}.
\end{equation}
\end{itemize}

\sset\subsubsection{}\label{subsec_descnt-data}
We shall call any pair $(X_\bullet,\omega_\bullet)$
of the above form, a {\em descent datum for the fibration $\phi$,
relative to the family $\underline S:=(\pi_i:S_i\to B~|~i\in I)$
and the cleavage $\sc$}.
The category of such descent data shall be denoted:
$$
\Desc(\phi,\underline S,\sc).
$$
Sometimes we may also denote it by $\Desc(\cA,\underline S,\sc)$,
if the notation is not ambiguous. Of course, two different choices
of cleavage lead to equivalent categories of descent data, so usually
we omit mentioning explicitly $\sc$, and write simply
$\Desc(\phi,\underline S)$ or $\Desc(\cA,\underline S)$.
The foregoing discussion can be summarized, by saying that
there is a commutative diagram of categories :
\set\begin{equation}\label{eq_replace-those-fctrs}
{\diagram
\cA(B) \ar[rrr]^-{\sCart_\cB(\iota_{\cS},\cA)}
\ar[d]_{\sev_B} & & & \sCart_\cB(\cS,\cA)
\ar[d]^-{\delta_{\underline S}} \\
\cA_B \ar[rrr]^-{\rho_{\underline S}} & & &
\Desc(\phi,\underline S,\sc)
\enddiagram}\end{equation}
whose vertical arrows are equivalences, and where
$\rho_{\underline S}$ is determined by $\sc$. Explicitly,
$\rho_{\underline S}$ assigns to every $C\in\Ob(\cA_B)$
the pair $(C_\bullet,\omega^C_\bullet)$ where $C_i:=\pi_i^*C$,
and $\omega^C_{ij}$ is the composition :
$$
\pi^{0*}_{ij}\circ\pi_j^* C
\xrightarrow[\sim]{\gamma_{(\pi^0_{ij},\pi_j)}}
(\pi_j\circ\pi^0_{ij})^*C=(\pi_i\circ\pi^1_{ij})^*C
\xrightarrow[\sim]{\gamma_{(\pi^1_{ij},\pi_i)}^{-1}}
\pi^{1*}_{ij}\circ\pi_i^* C
$$
where $\gamma_{(\pi^1_{ij},\pi_i)}$ and $\gamma_{(\pi^0_{ij},\pi_j)}$
are the coherence constraints for the cleavage $\sc$
(see \eqref{subsec_choice-of-pseudo}). The descent datum
$(X_\bullet,\omega_\bullet)$ is said to be {\em effective},
if it lies in the essential image of $\rho_{\underline S}$.

We also have an obvious functor :
$$
\sp_{\underline S}:\Desc(\phi,\underline S)\to
\prod_{i\in I}\cA_{S_i}
\qquad(X_i,\omega_{ij}~|~i,j\in I)\mapsto(X_i~|~i\in I)
$$
such that :
$$
\sp_{\underline S}\circ\rho_{\underline S}=
\prod_{i\in I}\pi_i^*:\cA_B\to\prod_{i\in I}\cA_{S_i}.
$$

\sset\subsubsection{}
Furthermore, for every morphism $f:B'\to B$ in $\cB$, set
$$
\underline S\times_Bf:=(\pi_i\times_BB':S_i\times_BB'\to B'~|~i\in I)
$$
which is a generating family for the sieve $\cS\times_Bf$
(notation of \eqref{subsec_slice-cat}); then we deduce a
pseudo-natural transformation of pseudo-functors
$(\cB/B)^o\to\bCat$ :
$$
\rho:\sc\circ i^o_B\Rightarrow\Desc(\phi,\underline S\times_B-,\sc)
\qquad
(f:B'\to B)\mapsto\rho_{\underline S\times_Bf}
$$
(where $\ss_B:\cB/B\to\cB$ is the source functor of
\eqref{subsec_slice-cat}) fitting into a commutative diagram :
$$
\xymatrix{ \cA(-)\circ\ss^o_B
\ar@{=>}[rrr]^-{\sCart_\cB(\iota_{\cS\times_B-},\cA)}
\ar@{=>}[d]_{\sev*\ss^o_B}
& & & \sCart_\cB(\cS\times_B-,\cA)
      \ar@{=>}[d]^{\delta_{\underline S\times_B-}} \\
\sc\circ\ss^o_B \ar@{=>}[rrr]^-\rho & & &
\Desc(\phi,\underline S\times_B-,\sc)
}$$
using which, one can figure out the pseudo-functoriality
of the rule : $f\mapsto\Desc(\phi,\underline S\times_Bf,\sc)$.
Namely, every pair of objects $f:C\to B$ and $f':C'\to B$,
and any morphism $h:C'\to C$ in $\cB/B$, yield a commutative
diagram :
$$
\xymatrix{ S_j\times_BC'
\ar[d]_{h_j:=S_j\times_Bh} & \ar[l]_-{\tilde\pi^0_{ij}}
S_{ij}\times_BC' \ar[r]^-{\tilde\pi^1_{ij}}
\ar[d]_{h_{ij}:=S_{ij}\times_Bh} &
S_i\times_BC' \ar[r]^-{\tilde\pi_i} \ar[d]^{h_i:=S_i\times_Bh} &
C' \ar[d]^h \\
S_j\times_BC & S_{ij}\times_BC \ar[l]^-{\pi^0_{ij}}
\ar[r]_-{\pi^1_{ij}} & S_i\times_BC \ar[r]_-{\pi_i} & C. }
$$
Hence one obtains a functor :
$$
\Desc(\phi,h,\sc):\Desc(\phi,\underline S\times_Bf,\sc)\to
                 \Desc(\phi,\underline S\times_Bf',\sc)
$$
by the rule :
$$
(X_i,\omega^X_{ij}~|~i,j\in I)\mapsto
(h_i^*X_i,\tilde\omega^X_{ij}~|~i,j\in I)
$$
where $\tilde\omega^X_{ij}$ is the isomorphism that makes
commute the diagram :
$$
\xymatrix{
h_{ij}^*\pi^{0*}_{ij}X_j
\ar[rr]^-{\gamma_{(h_{ij},\pi^0_{ij})}} \ar[d]_{h^*_{ij}\omega_{ij}} & &
(\pi^0_{ij}\circ h_{ij})^*X_j & & \tilde\pi^{0*}_{ij}\circ h_j^*X_j
\ar[ll]_-{\gamma_{(\tilde\pi^0_{ij},h_j)}} \ar[d]^{\tilde\omega^X_{ij}} \\
h_{ij}^*\pi^{1*}_{ij}X_i
\ar[rr]^-{\gamma_{(h_{ij},\pi^1_{ij})}} & &
(\pi^1_{ij}\circ h_{ij})^*X_i & &
\tilde\pi^{1*}_{ij}\circ h_i^*X_i
\ar[ll]_-{\gamma_{(\tilde\pi^1_{ij},h_i)}}
}$$
and if $f'':C''\to B$ is a third object, with a morphism
$g:C''\to C'$, we have a natural isomorphism of functors :
$$
\Desc(\phi,g,\sc)\circ\Desc(\phi,h,\sc)\Rightarrow\Desc(\phi,h\circ g,\sc)
$$
which is induced by the cleavage $\sc$, in the obvious fashion.

\begin{theorem}\label{th_first-theor}
For $i=1,2$, let $\phi_i:\cA_i\to\cB$ be two fibrations,
$F:\cA_1\to\cA_2$ a cartesian functor of $\cB$-categories,
$B$ an object of $\cB$ and $\cS$ a sieve of $\cB/B$ generated
by the family $(S_i\to B~|~i\in I)$. We assume that
$S_{ij}$ and $S_{ijk}$ are representable in $\cB$, for every
$i,j,k\in I$ (see remark {\em\ref{rem_abuse-not}}); then we have :
\begin{enumerate}
\item
For $n\in\{0,1,2\}$ and every $i,j,k\in I$, suppose that
\begin{enumerate}
\item
$\cS$ is a sieve both of\/ $\phi_1$-$n$-descent and
of\/ $\phi_2$-$(n-1)$-descent.
\item
The restriction $F_i:\phi_1^{-1}S_i\to\phi_2^{-1}S_i$ of $F$
is $n$-faithful.
\item
The restriction $F_{ij}:\phi_1^{-1}S_{ij}\to\phi_2^{-1}S_{ij}$
of $F$ is $(n-1)$-faithful.
\item
The restriction $F_{ijk}:\phi_1^{-1}S_{ijk}\to\phi_2^{-1}S_{ijk}$
of $F$ is $(n-2)$-faithful.
\end{enumerate}
Then the restriction $F_B:\phi^{-1}_1B\to\phi_2^{-1}B$ of $F$ is
$n$-faithful.
\item
Suppose that the functors $F_{ij}$ are fully faithful, and the
functors $F_{ijk}$ are faithful, for every $i,j,k\in I$. Then the
natural commutative diagram
$$
\xymatrix{
\sCart_\cB(\cS,\cA_1) \ar[rrr]^-{\sCart_\cB(\cS,F)} \ar[d] & & &
\sCart_\cB(\cS,\cA_2) \ar[d] \\
\prod_{i\in I}\phi_1^{-1}S_i \ar[rrr]^-{\prod_{i\in I}F_i} & & &
\prod_{i\in I}\phi_2^{-1}S_i 
}$$
is $2$-cartesian.
\end{enumerate}
\end{theorem}
\begin{proof}(i): In view of theorem \ref{th_split-fibration},
we may assume that both $\cA_1$ and $\cA_2$ are split fibrations
(with a suitable choice of cleavages), and $F$ is a split cartesian
functor. Recall that the latter condition means the following. For
every morphism $g:X\to Y$ in $\cB$, the induced diagram
$$
\xymatrix{
\phi^{-1}_1Y \ar[r]^-{g^*} \ar[d]_F & \phi^{-1}_1X \ar[d]^F \\
\phi^{-1}_2Y \ar[r]^-{g^*} & \phi^{-1}_2X 
}$$
commutes (where the horizontal arrows are the pull-back functor
given by the chosen cleavages). In this situation, we have a
commutative diagram
\set\begin{equation}\label{eq_reduce-to-desc}
{\diagram
\phi^{-1}_1B \ar[r]^-{\rho_{\underline S}} \ar[d]_{F_B} &
\Desc(\phi_1,\underline S) \ar[d]^{\underline F} \\
\phi^{-1}_2B \ar[r]^-{\rho_{\underline S}} &
\Desc(\phi_2,\underline S)
\enddiagram}
\end{equation}
whose right vertical arrow is the functor given by the rule :
$$
\underline X:=(X_i,\omega_{ij}~|~i,j\in I)\mapsto
\underline F(\underline X):=(F_iX_i,F_{ij}\omega_{ij}~|~i,j\in I).
$$
for every object $\underline X$ of $\Desc(\phi_1,\underline S)$.
By assumption, the top horizontal arrow of \eqref{eq_reduce-to-desc}
is $n$-faithful, and the bottom horizontal arrow is $(n-1)$-faithful.
We need to prove that the left vertical arrow is $n$-faithful,
and it is easily seen that this will follow, once we have shown
that the same holds for the right vertical arrow.

Suppose first that $n=0$; we have to show that
$\underline F$ is faithful. However, let $\underline X$ and
$\underline Y$ be two objects of $\Desc(\phi_1,\underline S)$,
and $\underline h_1,\underline h_2:\underline X\to\underline Y$
two morphisms. By definition, $\underline h_t$ (for $t=1,2$) is
a compatible system $(h_{t,i}:X_i\to Y_i~|~i\in I)$, where each
$h_{t,i}$ is a morphism in $\phi^{-1}_1S_i$. Then,
$\underline F(\underline h_t)$ is the compatible system
$(F_ih_{t,i}~|~i\in I)$. Thus, the condition
$\underline F(\underline h_1)=\underline F(\underline h_2)$
translates the system of identities $F_ih_{1,i}=F_ih_{2,i}$
for every $i\in I$. By assumption, each $F_i$ is faithful,
therefore $\underline h_1=\underline h_2$, as stated.

For $n=1$, assumption (d) is empty, (b) means that $F_i$ is
fully faithful, and (c) means that $F_{ij}$ is faithful for
every $i,j\in I$. In light of the previous case, we have only
to show that $\underline F$ is full. Hence, let
$\underline X,\underline Y$ be as in the foregoing, and
$(h_i:F_iX_i\to F_iY_i~|~i\in I)$ a morphism
$\underline F(\underline X)\to\underline F(\underline Y)$ in
$\Desc(\phi_2,\underline S)$. By assumption, for every $i\in I$
we may find a unique morphism $f_i:X_i\to Y_i$ such that
$F_if_i=h_i$. It remains only to check that the system
$(f_i~|~i\in I)$ fulfills condition \eqref{eq_morphius},
and since the functors $F_{ij}$ are faithful, it suffices
to verify that $F_{ij}\eqref{eq_morphius}$ holds. However,
since $F$ is split cartesian, we have :
$$
F_{ij}\circ\pi^{0*}_{ij}f_j=\pi^{0*}_{ij}\circ F_jf_j
\qquad
F_{ij}\circ\pi^{1*}_{ij}f_i=\pi^{1*}_{ij}\circ F_if_i
$$
hence we reduce to showing that
$(F_{ij}\omega^Y_{ij})\circ\pi^{0*}_{ij}h_j=
\pi^{1*}_{ij}h_i\circ(F_{ij}\omega^X_{ij})$, which holds by
assumption.

Next, we consider assertion (ii) : the contention is that the
functors $\underline F$ and :
$$
\sp_{1,\underline S}:\Desc(\phi_1,\underline S)\to
\prod_{i\in I}\phi_1^{-1}S_i
$$
as in \eqref{subsec_descnt-data}, induce an equivalence
$(\sp_{1,\underline S},\underline F)$ between
$\Desc(\phi_1,\underline S)$ and the category $\cC$ consisting
of all data of the form
$\underline G:=(G_i,H_i,\alpha_i,\omega^H_{ij}~|~i,j\in I)$,
where $G_i\in\Ob(\phi_1^{-1}S_i)$, $H_i\in\Ob(\phi^{-1}_2S_i)$,
$\alpha_i:F_i G_i\isom H_i$ are isomorphisms in $\phi_2^{-1}S_i$,
and $\underline H:=(H_i,\omega^H_{ij}~|~i,j\in I)$ is an object
of $\Desc(\phi_2,\underline S)$. However, given an object as above,
set :
$$
\omega^{H'}_{ij}:=(\pi^{1*}_{ij}\alpha_i^{-1})\circ
\omega^H_{ij}\circ(\pi^{0*}_{ij}\alpha_j).
$$
Since $\phi_2$ is a split fibration, one verifies easily that the
datum $\underline H':=(H'_i:=F_iG_i,\omega^{H'}_{ij}~|~i,j\in I)$
is an object of $\Desc(\phi_2,\underline S)$ isomorphic to
$\underline H$, and the new datum
$(F_i,H'_i,\one_{H'_i},\omega^{H'}_{ij}~|~i,j\in I)$
is isomorphic to $\underline G$; hence $\cC$ is equivalent
to the category $\cC'$ whose objects are all data of the form
$(G_i,\omega_{ij}~|~i,j\in I)$ where $G_i\in\Ob(\phi_1^{-1}S_i)$,
and $(F_iG_i,\omega_{ij}~|~i,j\in I)$ is an object of
$\Desc(\phi_2,\underline S)$. By assumption $F_{ij}$ is fully
faithful, and $F$ is a split cartesian functor; hence we may
find unique isomorphisms
$\tilde\omega^G_{ij}:\pi^{0*}_{ij}G_j\isom\pi^{1*}_{ij}G_i$
such that $\tilde\omega_{ij}=F_{ij}\tilde\omega^G_{ij}$. We
claim that the datum $(G_i,\omega^G_{ij}~|~i,j\in I)$
is an object of $\Desc(\phi_1,\underline S)$, {\em i.e.}
the isomorphisms $\omega^G_{ij}$ satisfy the cocycle condition
\set\begin{equation}\label{eq_new-cocycle}
\pi^{2*}_{ijk}\omega^G_{ij}\circ\pi^{0*}_{ijk}\omega^G_{jk}=
\pi^{1*}_{ijk}\omega^G_{ik}
\qquad
\text{for every $i,j,k\in I$}.
\end{equation}
To check this identity, since by assumption the
functors $F_{ijk}$ are faithful, it suffices to see that
$F_{ijk}\eqref{eq_new-cocycle}$ holds, which is clear, since the
cocycle condition holds for the isomorphisms $\omega_{ij}$
(and since $F$ is split cartesian).
This shows that $(\sp_{1,\underline S},\underline F)$ is essentially
surjective. Next, since the functor $\sp_{1,\underline S}$ is
faithful, the same holds for $(\sp_{1,\underline S},\underline F)$.
Finally, let
$$
\underline G:=(G_i,\omega_{ij}~|~i,j\in I)
\qquad
\underline G':=(G'_i,\omega'_{ij}~|~i,j\in I)
$$
be two objects of $\cC'$. A morphism $\underline G\to\underline G'$
consists of a system $(\alpha_i:G_i\to G'_i~|~i\in I)$
of morphisms such that $(F_i\alpha_i~|~i\in I)$ is a morphism
$$
(F_iG_i,\omega_{ij}~|~i,j\in I)\to(F_iG'_i,\omega'_{ij}~|~i,j\in I)
$$
in $\Desc(\phi_2,\underline S)$. To show that
$(\sp_{1\underline S},\underline F)$ is full, and since we know
already that this functor is essentially surjective, we may
assume that there exist $(G_i,\omega^G_{ij}~|~i,j\in I)$,
$(G'_i,\omega^{G'}_{ij}~|~i,j\in I)$ in $\Desc(\phi_1,\underline S)$
such that $\omega_{ij}=F_{ij}\omega^G_{ij}$ and
$\omega'_{ij}=F_{ij}\omega^{G'}_{ij}$ for every $i,j\in I$; in this
case, it suffices to verify the identity
\set\begin{equation}\label{eq_alphamorphius}
\omega^{G'}_{ij}\circ\pi^{0*}_{ij}\alpha_j=
\pi^{1*}_{ij}\alpha_i\circ\omega^G_{ij}
\qquad
\text{for every $i,j\in I$}.
\end{equation}
Again, the faithfulness of $F_{ij}$ reduces to checking that
$F_{ij}\eqref{eq_alphamorphius}$ holds, which is clear, since $F$
is split cartesian.

Lastly, notice that the case $n=2$ of assertion (i) is a
formal consequence of (ii).
\end{proof}

\sset\subsubsection{}\label{subsec_th-descent}
In the situation of \eqref{subsec_descnt-data}, let $\cS$
be the sieve generated by the family $\underline S$, and
$g:B'\to B$ any morphism in $\cB$. We let :
$$
B'_i:=B'\times_BS_i
\quad
B'_{ij}:=B'\times_BS_{ij}
\quad
B'_{ijk}:=B'\times_BS_{ijk}
\qquad
\text{for every $i,j,k\in I$}
$$
and denote $g_i:B'_i\to S_i$, $g_{ij}:B'_{ij}\to S_{ij}$ and
$g_{ijk}:B'_{ijk}\to S_{ijk}$ the induced projections.

\begin{corollary}\label{cor_first-theor}
With the notation of \eqref{subsec_th-descent}, let $n\in\{0,1,2\}$.
The following holds :
\begin{enumerate}
\item
$\cS$ is a sieve of $\phi$-$n$-descent, if and only if
$\rho_{\underline S}$ is $n$-faithful (see
\eqref{eq_replace-those-fctrs}). 
\item
Suppose that :
\begin{enumerate}
\item
$\cS$ is a sieve of universal $\phi$-$n$-descent.
\item
The pull-back functors $g^*_i:\phi^{-1}S_i\to\phi^{-1}B'_i$
are $n$-faithful.
\item
The pull-back functors
$g^*_{ij}:\phi^{-1}S_{ij}\to\phi^{-1}B'_{ij}$ are $(n-1)$-faithful.
\item
The pull-back functors
$g^*_{ijk}:\phi^{-1}S_{ijk}\to\phi^{-1}B'_{ijk}$ are $(n-2)$-faithful.
\end{enumerate}
Then the pull-back functor $g^*$ is $n$-faithful.
\item
Suppose that the functors $g_{ij}^*$ are fully faithful,
and the functors $g_{ijk}^*$ are faithful, for every $i,j,k\in I$.
Then the natural essentially commutative diagram :
$$
\xymatrix{\Desc(\phi,\underline S)
\ar[rr]^-{\Desc(\phi,g)} \ar[d]_{\sp_{\underline S}} & &
\Desc(\phi,\underline S\times_BB') \ar[d]^{\sp_{\underline S\times_BB'}} \\
\prod_{i\in I}\phi^{-1}S_i \ar[rr]^-{\prod_{i\in I}g^*_i} & &
\prod_{i\in I}\phi^{-1}B'_i
}$$
is $2$-cartesian (see remark {\em\ref{rem_pseudo-limit}(v)}
and example {\em\ref{ex_2-products}(ii)}).
\end{enumerate}
\end{corollary}
\begin{proof} (i) follows by inspecting \eqref{eq_replace-those-fctrs}.

(ii): Thanks to theorem \ref{th_split-fibration}, we may assume
that $\phi$ is a split fibration. Now, set
$$
\cC:=\sMorph(\cB)
\qquad
\cA_1:=\cC\times_{(\st,\phi)}\cA
\qquad
\cA_2:=\cC\times_{(\ss,\phi)}\cA
$$
where $\ss,\st:\cC\to\cB$ are the source and target functors
(see \eqref{subsec_Morph-cat}). The natural projections
$\phi_i:\cA_i\to\cC$ (for $i=1,2$) are two fibrations (see
remark \ref{rem_added-little-extra}(i)). Moreover, $\st$
induces a functor $\st_{|g}:\cC\!/g\to\cB\!/\!B$ (notation
of \eqref{eq_restrict-over-X}) and we let
$\cS\!/g:=\st_{|g}^{-1}\cS$, which is the sieve of $\cC\!/g$
generated by the cartesian diagrams
$$
D_i \quad : \quad
{\diagram B'_i \ar[r]^-{g_i} \ar[d] & S_i \ar[d] \\
          B' \ar[r]^-g & B \enddiagram}
\qquad
\text{for every $i\in I$}.
$$
Notice that the products $D_{ij}:=D_i\times D_j$ are
represented in $\cC\!/g$ by the diagrams
$$
{\diagram B'_{ij} \ar[r]^-{g_{ij}} \ar[d] & S_{ij} \ar[d] \\
          B' \ar[r]^-g & B \enddiagram}
\qquad
\text{for every $i,j\in I$}
$$
and likewise one may represent the triple products
$D_{ijk}:=D_{ij}\times D_k$.

By definition, the objects of $\cA_1$ (resp. $\cA_2$) are the
pairs $(h:X\to Y,a)$, where $h$ is a morphism in $\cB$ and
$a\in\Ob(\phi^{-1}Y)$ (resp. $a\in\Ob(\phi^{-1}X)$). A morphism
of $\cA_1$ (resp. of $\cA_2$)
$$
(h:X\to Y,a)\to(h':X'\to Y',a')
$$
is a datum $(f_1,f_2,t)$, where $f_1:X\to X'$ and $f_2:Y\to Y'$
are morphisms in $\cB$ with $f_2\circ h=h'\circ f_1$, and
$t:a\to f_2^*a'$ (resp. $t:a\to f_1^*a'$) is a morphism in
$\phi^{-1}Y$ (resp. in $\phi^{-1}X$). Now, we define a functor
$F:\cA_1\to\cA_2$ of $\cC$-categories, by the rule :
\begin{itemize}
\item
$(h,a)\mapsto(h,h^*a)$ for every $(h,a)\in\Ob(\cA_1)$.
\item
$(f_1,f_2,t)\mapsto(f_1,f_2,h^*t)$ for every morphism
$(f_1,f_2,t)$ of $\cA_1$ as above. Notice that if
$t:a\to f_2^*a'$ is a morphism in $\phi^{-1}Y$, then
$h^*t:h^*a\to h^*f^*_2a'=f_1^*h'{}^*a'$ is a morphism
of $\phi^{-1}Y'$, since $\phi$ is a split fibration.
\end{itemize}
Notice that a morphism $(f_1,f_2,t)$ of either $\cA_1$
or $\cA_2$ is cartesian if and only if $t$ is an
isomorphism; especially, it is clear that $F$ is a
cartesian functor. Moreover, for every object $h:X\to Y$
of $\cC$, the restriction $\phi_1^{-1}h\to\phi^{-1}_2h$
of $F$ is isomorphic to the pull-back functor
$h^*:\phi^{-1}Y\to\phi^{-1}h$. Especially, conditions
(b)--(d) say that the restriction
$F_i:\phi^{-1}_1g_i\to\phi^{-1}_2g_i$ (resp.
$F_{ij}:\phi^{-1}_1g_{ij}\to\phi^{-1}_2g_{ij}$, resp.
$F_{ijk}:\phi^{-1}_1g_{ijk}\to\phi^{-1}_2g_{ijk}$) are
$n$-faithful (resp. $(n-1)$-faithful, resp. $(n-2)$-faithful).
In light of theorem \ref{th_first-theor}(i), we are then
reduced to showing

\begin{claim}\label{cl_both-desc}
$\cS/g$ is a sieve both of $\phi_1$-$n$-descent and of
$\phi_2$-$n$-descent.
\end{claim}
\begin{pfclaim}[] Let $\cD$ be any (small) category; we remark
first that a functor $\cD\to\cA_1$ is the same as a pair of
functors $(H:\cD\to\cA,K:\cD\to\cC)$ such that
$\phi\circ H=\st\circ K$, and likewise one can describe the
functors $\cD\to\cA_2$. Then, it is easily seen that the
functors
$$
\begin{aligned}
\cA(B)\to\sCart_\cC(\cC/g,\cA_1) \qquad &
G\mapsto(G\circ\st,\st) \\
\cA(B')\to\sCart_\cC(\cC/g,\cA_2) \qquad &
G\mapsto(G\circ\ss,\ss)
\end{aligned}
$$
are equivalences, and induce equivalences
$$
\begin{aligned}
\sCart_\cB(\cS,\cA)\to\sCart_\cC(\cS/g,\cA_1) \\
\sCart_\cB(\cS\times_Bg,\cA)\to\sCart_\cC(\cS/g,\cA_2)
\end{aligned}
$$
(details left to the reader). The claim follows immediately.
\end{pfclaim}
\end{proof}

\sset\subsubsection{}\label{subsec_fibrat-desc-cat}
Quite generally, if $\phi:\cA\to\cB$ is a fibration over a
category $\cB$ that admits fibre products, the descent data
for $\phi$ (relative to a fixed cleavage $\sc$) also form a
fibration :
$$
\sD\phi:\phi\tdu\Desc\to\sMorph(\cB).
$$
Namely, for every morphism $f:T'\to T$ of $\cB$, the fibre
over $f$ is the category $\Desc(\phi,f)$ of all descent data
$(f,A,\xi)$ relative to the family $\{f\}$, and the cleavage
$\sc$, so $A$ is an object of $\phi^{-1}T'$ and
$\xi:p^*_1 A\isom p^*_2 A$ is an isomorphism in the category
$\phi^{-1}(T'\times_TT')$ satisfying the usual cocycle condition
(here $p_1,p_2:T'\times_TT'\to T'$ denote the two natural
morphisms). Given two objects $\underline A:=(f:T'\to T,A,\xi)$,
$\underline A':=(g:W'\to W,A',\zeta)$ of $\phi\tdu\Desc$, the
morphisms $\underline A\to\underline A'$ are the data $(h,\alpha)$
consisting of a commutative diagram :
$$
\xymatrix{
W' \ar[r]^h \ar[d]_g & T' \ar[d]^f \\ W \ar[r] & T
}$$
and a morphism $\alpha:A\to A'$ such that $\phi(\alpha)=h$
and $p^*_1(\alpha)\circ\xi=\zeta\circ p_2^*(\alpha)$.

We have a natural cartesian functor of fibrations :
$$
\xymatrix{
\cA\times_{(\phi,\st)}\sMorph(\cB) \ar[rr]^-\sd \ar[rd]_\sp & &
\phi\tdu\Desc \ar[dl]^{\sD\phi} \\
& \sMorph(\cB)
}$$
where $\st:\sMorph(\cB)\to\cB$ is the target functor (see
\eqref{subsec_Morph-cat} and example \ref{ex_cat-cats}(i)).
Namely, to any pair $(T,f:S\to\phi T)$ with $T\in\Ob(\cA)$
and $f\in\Ob(\sMorph(\cB))$, one assigns the canonical
descent datum $\sd(T,f):=\rho_{\{f\}}(T)$ in $\Desc(\phi,f)$
associated with the pair as in \eqref{eq_replace-those-fctrs}.

\begin{corollary}\label{cor_descend-descent}
In the situation of \eqref{subsec_fibrat-desc-cat}, let
$f:B'\to B$ be a morphism of $\cB$, and $\cS$ a sieve
of $\cB/B$, generated by a family $(S_i\to B~|~i\in I)$.
Let $n\in\{0,1,2\}$, and suppose that :
\begin{enumerate}
\alphaenu
\item
$\cS$ is a sieve of universal $\phi$-$n$-descent.
\item
For every $i\in I$, the morphism $S_i\times_Bf$ is of
$\phi$-$n$-descent.
\item
For every $i,j\in I$, the morphism $S_{ij}\times_Bf$ is of
$\phi$-$(n-1)$-descent.
\item
For every $i,j,k\in I$, the morphism $S_{ijk}\times_Bf$ is
of $\phi$-$(n-2)$-descent.
\end{enumerate}
Then $f$ is a morphism of $\phi$-$n$-descent.
\end{corollary}
\begin{proof} In view of corollary \ref{cor_first-theor}(i),
it is easily seen that a morphism $g:T'\to T$ in $\cB$ is of
$\phi$-$n$-descent if and only if the restriction
$\phi^{-1}T\to\sD\phi^{-1}g$ of $\sd$ is $n$-faithful.
Set $\cC:=\sMorph(\cB)$; as in the proof of corollary
\ref{cor_first-theor}(ii), the functor $\st$ induces a
functor $\st_{|f}:\cC/f\to\cB/B$, and we let
$\cS\!/\!f:=\st_{|f}^{-1}\cS$. With this notation,
theorem \ref{th_first-theor}(i) reduces to showing :
\begin{claim} The sieve $\cS\!/\!f$ is both of $\sp$-$n$-descent
and of $\sD\phi$-$n$-descent.
\end{claim}
\begin{pfclaim}[] By claim \ref{cl_both-desc}, it is already known
that $\cS/f$ is of $\sp$-$n$-descent. To show that $\cS/f$
is of $\sD\phi$-$n$-descent, we consider the commutative diagram
$$
\xymatrix{
\cA(B') \ar[r] \ar[d] &
\sCart_\cC(\cC/f,\phi\tdu\Desc) \ar[d] \\
\sCart_\cB(\cS\times_Bf,\cA) \ar[r] &
\sCart_\cC(\cS/f,\phi\tdu\Desc)
}$$
whose left (resp. right) vertical arrow is induced by the inclusion
$\cS\times_Bf\to\cB/B'$ (resp. $\cS/f\to\cS\times_Bf$) and whose
top horizontal arrow is defined as follows. Given a cartesian
functor $G:\cB/B'\to\cA$, we let $\sD G:\cC/f\to\phi\tdu\Desc$
be the unique cartesian functor determined on the objects of
$\cC\!/\!f$ by the rule :
$$
\left({\diagram T' \ar[r]^-g \ar[d]_h & T \ar[d] \\
           B' \ar[r]^-f & B
\enddiagram}\right)
\mapsto\sd(G(h),g).
$$
We leave to the reader the verification the rule $G\mapsto\sD G$
extends to a well defined functor, and then there exists a unique
(similarly defined) bottom horizontal arrow that makes commute the
foregoing diagram. Moreover, both horizontal arrow thus obtained
are equivalences of categories. The claim follows.
\end{pfclaim} 
\end{proof}

\begin{lemma}\label{lem_was-pair-of-sieves}
Let $i\leq 2$ be an integer, $F:\cE\to\cC$ a fibration,
$X\in\Ob(\cC)$, and $\cT\subset\cS$ two sieves of\/ $\cC/X$.
The following holds :
\begin{enumerate}
\item
If $\cT\times_Xf$ is of $F$-$i$-descent for every
$(Y\xrightarrow{f}X)\in\cS$, then the induced functor
$$
j:\sCart_\cC(\cS,\cE)\to\sCart_\cC(\cT,\cE)
$$
is $i$-faithful.
\item
If $\cT$ is of universal $F$-$i$-descent, the same holds for $\cS$.
\end{enumerate}
\end{lemma}
\begin{proof}(i): The assertion is trivial for $i<0$. We consider
first the case where $i=0$; thus, let $\phi,\psi:\cS\to\cE$
be two $\cC$-cartesian functors, and
$\alpha,\beta:\phi\Rightarrow\psi$ two natural
$\cC$-transformations such that $\alpha_{|\cT}=\beta_{|\cT}$.
Let also $(Y\xrightarrow{f}X)\in\Ob(\cS)$ be any object
and $f_*:\cC/Y\to\cS$ the induced functor; then
$\cT':=\cT\times_Xf\subset\cC/Y$ is a sieve of
$F$-$0$-descent and
$$
(\alpha*f_*)_{|\cT'}=(\beta*f_*)_{|\cT'}.
$$
We deduce that $\alpha*f_*=\beta*f_*$; especially
$\alpha_f=\beta_f$, which shows that $j$ is faithful.

Next, suppose that $i=1$, and for $\phi$ and $\psi$ as in the
foregoing, let $\alpha:\phi_{|\cT}\Rightarrow\psi_{|\cT}$ be a
given natural $\cC$-transformation, and
$(Y\xrightarrow{f}X)\in\Ob(\cS)$; then the sieve
$\cT':=\cT\times_Xf\subset\cC/Y$ is of $F$-$1$-descent,
and $f_*$ restricts to a functor $f_{*|\cT'}:\cT'\to\cT$,
so the natural $\cC$-transformation
$$
\alpha*(f_{*|\cT'}):
(\phi\circ f_*)_{|\cT'}\Rightarrow(\psi\circ f_*)_{|\cT'}
$$
extends to a unique natural $\cC$-transformation
$\alpha^{(f)}:\phi\circ f_*\Rightarrow\psi\circ f_*$.
Now, let us define
$$
\tilde\alpha_f:=\alpha^{(f)}_{\one_Y}
\qquad
\text{for every $(Y\xrightarrow{f}X)\in\Ob(\cS)$}.
$$
Clearly, $\tilde\alpha_f=\alpha_f$ whenever $f\in\Ob(\cT)$,
and it remains to check that $\tilde\alpha$ is a natural
$\cC$-transformation $\phi\Rightarrow\psi$. Thus, let
$h/X:(Y'\xrightarrow{f'}X)\to(Y\xrightarrow{f}X)$ be a
morphism of $\cS$; we need to show that
$$
\psi(h/X)\circ\tilde\alpha_f=\tilde\alpha_{f'}\circ\phi(h/X).
$$
Considering the morphism $h/Y:(Y'\xrightarrow{h}Y)\to\one_Y$
of $\cC/Y$, and noticing that $\phi(h/X)=(\phi\circ f_*)(h/Y)$
and $\psi(h/X)=(\psi\circ f_*)(h/Y)$, we come down to checking
that $\alpha^{(f)}_h=\alpha^{(f')}_{\one_{Y'}}$. We show more
precisely that $\alpha^{(f)}*h_*=\alpha^{(f')}$. Indeed, notice that
$$
(\alpha^{(f)}*h_*)_{|\cT\times_Xf'}=\alpha*(f'_{*|\cT\times_Xf'})=
\alpha^{(f')}_{|\cT\times_Xf'}
$$
and since $\cT\times_Xf'$ is of $F$-0-descent, the assertion
follows.

Lastly, for $i=2$, we know already that $j$ is fully faithful,
so it remains only to show that every cartesian functor
$\phi:\cT\to\cE$ is isomorphic to the restriction of a
cartesian functor $\cS\to\cE$. Now, for every
$(f:Y\to X)\in\Ob(\cS)$ set $\cT^{(f)}:=\cT\times_Xf$; by
assumption there exists a $\cC$-cartesian functor
$\phi^{(f)}:\cC/Y\to\cE$ with an isomorphism of functors :
$$
\omega_f:\phi^{(f)}_{|\cT^{(f)}}\isom\phi\circ(f_{*|\cT^{(f)}}).
$$
Notice that if $f\in\Ob(\cT)$, we have $\cT^{(f)}=\cC/Y$,
and in this case we take $\phi^{(f)}:=\phi\circ f_*$ and
$\omega_f:=\one_{\phi^{(f)}}$. Then, for every morphism
$h/X:(Y'\xrightarrow{f'}X)\to(Y\xrightarrow{f'}X)$ in $\cS$
we get the isomorphism :
$$
\phi^{(f')}_{|\cT^{(f')}}\xrightarrow{\omega_{f'}}
\phi\circ f'_{*|\cT^{(f')}}\xrightarrow{(\omega_f*h_*)_{|\cT^{(f')}}^{-1}}
(\phi^{(f)}\circ h_*)_{*|\cT^{(f')}}
$$
which, as $\cT^{(f')}$ is of $F$-1-descent, extends uniquely
to an isomorphism of functors $\cC/Y'\to\cE$
$$
\tau_{h/X}:\phi^{(f')}\isom\phi^{(f)}\circ h_*.
$$

\begin{claim}\label{cl_snow}
The rule : $(f:Y\to X)\mapsto(\phi^{(f)}:\cC/Y\to\cE)$
for every $f\in\Ob(\cS)$ defines a pseudo-cocone
$$
\phi^{(\bullet)}:G_\cS\Rightarrow\sF_\cE
$$
whose coherence constraint is given by the system of isomorphisms
$\tau_{\bullet/X}$, and where $G_\cS:\cS\to\bCat$ is the functor
associated with $\cS$ as in \eqref{subsec_towards-top}.
\end{claim}
\begin{pfclaim} Clearly
$(\tau_{\one_Y/Y})_{|\cT^{(f)}}=\one_{\phi^{(f)}_{|\cT^{(f)}}}$ for every
$(f:Y\to X)\in\Ob(\cS)$, whence $\tau_{\one_Y/Y}=\one_{\phi^{(f)}}$
by the uniqueness of $\tau_{\one_Y/Y}$. Next, let
$(Y''\xrightarrow{f''}X)\xrightarrow{h'/X}(Y'\xrightarrow{f'}X)
\xrightarrow{h/X}(Y\xrightarrow{f}X)$ be two morphism of $\cS$;
the coherence axiom for $\tau_{\bullet/X}$ comes down to the
identity :
$$
(\tau_{h/X}*h'_*)\circ\tau_{h'/X}=\tau_{(h\circ h')/X}.
$$
But the latter can be checked after restriction to the sieve
$\cT^{(f'')}$, where it follows by a simple inspection : the
details are left to the reader.
\end{pfclaim}

From claim \ref{cl_snow} it follows that there exists a
functor
$$
\tilde\phi:\cS\to\cE
\qquad\text{and an invertible modification}\qquad
\Xi:\sF_{\tilde\phi}\odot\hat\eps\isom\phi^{(\bullet)}
$$
where $\hat\eps:G_\cS\Rightarrow\sF_\cE$ is the universal
pseudo-cocone provided by lemma \ref{lem_new-pseudo-col}.
Explicitly, $\Xi$ is a system of isomorphisms of functors
$$
\Xi_f:\tilde\phi\circ f_*\isom\phi^{(f)}
\qquad
\text{for every $(f:Y\to X)\in\Ob(\cS)$}
$$
from which it follows easily that $\tilde\phi$ is
$\cC$-cartesian. The compatibility conditions for $\Xi$
amount to the identities :
\set\begin{equation}\label{eq_compat-Xi}
\Xi_f*h_*=\tau_h\odot\Xi_{f'}
\qquad
\text{for every morphism $h/X:f\to f'$ in $\cS$}.
\end{equation}
To conclude, it suffices to exhibit an isomorphism of
functors $\xi:\tilde\phi_{|\cT}\isom\phi$. To this aim, we set
$$
\xi_g:=(\Xi_g)_{\one_Y}:\tilde\phi(g)\isom\phi^{(g)}(\one_Y)=\phi(g)
\qquad
\text{for every $(g:Y\to X)\in\Ob(\cT)$}.
$$
To check the naturality of $\xi$, let
$h/X:(Y'\xrightarrow{g'}X)\to(Y\xrightarrow{g}X)$ be any
morphism of $\cT$; notice that $\tau_{h/X}=\one_{\phi^{(g')}}$,
whence $\xi_{g'}=(\Xi_g)_h$, due to \eqref{eq_compat-Xi}.
On the other hand, the naturality of $\Xi_g$ yields the
commutative diagram
$$
\xymatrix{
\tilde\phi(g')=\tilde\phi\circ g_*(h) \ar[rr]^-{(\Xi_g)_h}
\ar[d]_{\tilde\phi\circ g_*(h/Y)} & & \phi^{(g)}(h)=\phi(g')
\ar[d]^{\phi^{(g)}(h/Y)} \\
\tilde\phi(g)=\tilde\phi\circ g_*(\one_Y) \ar[rr]^-{\xi_g} & &
\phi^{(g)}(\one_Y)=\phi(g)
}$$
where $\tilde\phi\circ g_*(h/Y)=\tilde\phi(h/x)$ and
$\phi^{(g)}(h/Y)=\phi(h/X)$, whence the contention.

(ii): We consider the induced functors
$$
\cE(X)\xrightarrow{j'}\sCart_\cC(\cS,\cE)
\xrightarrow{j''}\sCart_\cC(\cT,\cE)
$$
and notice that $j''$ is $i$-faithful, by virtue of (i);
the same holds for $j''\circ j'$, by assumption. Then it
follows easily that $j'$ is $i$-faithful as well, whence
the assertion.
\end{proof}

\begin{proposition}\label{prop_topol-of-univ-descent}
Let $i\leq 2$ be an integer, $F:\cE\to\cC$ a fibration,
$X\in\Ob(\cC)$, and $\cT,\cS$ two sieves of\/ $\cC/X$,
such that :
\begin{enumerate}
\alphaenu
\item
$\cS$ is of universal $F$-i-descent.
\item
$\cT\times_Xf$ is of universal $F$-$i$-descent for every
$(Y\xrightarrow{f}X)\in\cS$.
\end{enumerate}
Then $\cT$ is of universal $F$-$i$-descent.
\end{proposition}
\begin{proof} According to lemma \ref{lem_was-pair-of-sieves}(ii),
the sieve $\cS':=\cS\bigcup\cT$ is also of universal $F$-$i$-descent.
It is also clear that $\cT\times_Xf$ is of $F$-$i$-descent for
every $(Y\xrightarrow{f}X)\in\cS'$. Thus, we may replace $\cS$
by $\cS'$ and assume that $\cT\subset\cS$. Next, let $f:Y\to X$
be any morphism of $\cC$; we need to check that $\cT\times_Xf$
is a sieve of $F$-$i$-descent, and by assumption the same holds
for $\cS\times_Xf$, so it suffices to show that the induced
functor
$\sCart_\cC(\cS\times_Xf,\cE)\to\sCart_\cC(\cT\times_Xf,\cE)$ is
$i$-faithful. In view of lemma \ref{lem_was-pair-of-sieves}(i),
we are then reduced to checking that $(\cT\times_Xf)\times_Yg$
is a sieve of $F$-$i$-descent for every $g\in\Ob(\cS\times_Xf)$.
But we have $(\cT\times_Xf)\times_Yg=\cT\times_X(f\circ g)$,
and $f\circ g\in\Ob(\cS)$, so the assertion follows from (b).
\end{proof}

\subsection{Profinite groups and Galois categories}
\label{sec_profinite-gens}
Quite generally, for any profinite group $P$, let $P\tdu\Set$
denote the category of discrete finite sets, endowed with a
continuous left action of $P$ (the morphisms in $P\tdu\Set$
are the $P$-equivariant maps). Any continuous group homomorphism
$\omega:P\to Q$ of profinite groups induces a {\em restriction
functor}
$$
\Res(\omega):Q\tdu\Set\to P\tdu\Set
$$
in the obvious way. In case the notation is not ambiguous,
one writes also $\Res^P_Q$ for this functor.

For any two profinite groups $P$ and $Q$, we denote by
$$
\Hom_\mathrm{cont}(P,Q)
$$
the set of all continuous group homomorphisms $P\to Q$. If
$\phi_1,\phi_2$ are two such group homomorphisms, we say that
$\phi_1$ is {\em conjugate\/} to $\phi_2$, and we write
$\phi_1\sim\phi_2$, if there exists an inner automorphism
$\omega$ of $G$, such that $\phi_2=\omega\circ\phi_1$. Clearly
the trivial map $\pi\to G$ (whose image is the neutral element
of $G$), is the unique element of a distinguished conjugacy class.

\sset\subsubsection{}\label{subsec_first-non-abel}
Let $P$ be any profinite group; for any (discrete) finite group
$G$, consider the pointed set $\Hom_\mathrm{cont}(P,G)/\!\!\sim$
of conjugacy classes of continuous group homomorphisms $P\to G$.
This is also denoted
$$
H^1_\mathrm{cont}(P,G)
$$
and called the first {\em non-abelian continuous cohomology group\/}
of $P$ with coefficients in $G$ (so $G$ is regarded as a $P$-module
with trivial $P$-action). Clearly the formation of $H^1(P,G)$ is
covariant on the argument $G$, and controvariant for continuous
homomorphisms of profinite groups.

\begin{lemma}\label{lem_profinite-iso-crit}
Let $\phi:P\to P'$ be a continuous homomorphism of profinite groups,
and suppose that the induced map of pointed sets :
$$
H^1_\mathrm{cont}(P',G)\to H^1_\mathrm{cont}(P,G)
\quad :\quad f\mapsto f\circ\phi
$$
is bijective, for every finite group $G$. Then $\phi$ is an isomorphism
of topological groups.
\end{lemma}
\begin{proof} First we show that $\phi$ is injective. Indeed,
let $x\in P$ be any element; we may find an open normal subgroup
$H\subset P$ such that $x\notin H$; taking $G:=P/H$, we deduce that
the projection $P\to P/H$ factors through $\phi$ and a group
homomorphism $f:P'\to P/H$, hence $x\notin\Ker\,\phi$, as claimed.
Moreover, let $H':=\Ker\,f$; clearly $H'\cap P=H$, so the topology
of $P$ is induced from that of $P'$. It remains only to show that
$\phi$ is surjective; to this aim, we consider any continuous
surjection $f':P'\to G'$ onto a finite group, and it suffices to
show that the restriction of $f'$ to $\phi P$ is still surjective.
Indeed, let $G$ be the image of $\phi P$ in $G'$, denote by
$i:G\to G'$ the inclusion map, and let $f:P\to G$ be the unique
continuous map such that $i\circ f=f'\circ\phi$; by assumption,
there exists a continuous group homomorphism $g:P'\to G$ such that
$f=g\circ\phi$. On the other hand, $(i\circ g)\circ\phi=f'\circ\phi$,
hence the conjugacy class of $i\circ g$ equals the conjugacy class
of $f'$, especially $i\circ g$ is surjective, hence the same holds
for $i$, as required.
\end{proof}

\sset\subsubsection{}
Let $G$ be a profinite group, and $H\subset G$ an open subgroup.
It is easily seen that $H$ is also a profinite group, with the
topology induced from $G$. Morever, the restriction functor
$$
\Res^H_G:G\tdu\Set\to H\tdu\Set
$$
admits a left adjoint
$$
\Ind^G_H:H\tdu\Set\to G\tdu\Set.
$$
Namely, to any finite set $\Sigma$ with a continuous left action
of $H$, one assigns the set $\Ind^G_H\Sigma:=G\times\Sigma/\!\sim$,
where $\sim$ is the equivalence relation such that
$$
(gh,\sigma)\sim(g,h\sigma)
\qquad
\text{for every $g\in G$, $h\in H$ and $\sigma\in\Sigma$.}
$$
The left $G$-action on $\Ind^G_H\Sigma$ is given by the rule :
$(g',(g,\sigma))\mapsto(g'g,\sigma)$ for every $g,g'\in G$ and
$\sigma\in\Sigma$. It is easily seen that this action is continuous,
and the reader may check that the functor $\Ind^G_H$ is indeed
left adjoint to $\Res^H_G$.

\sset\subsubsection{}\label{subsec_nat-quinverse}
Moreover, let $1$ denote the final object of $H\tdu\Set$;
notice that $\Ind^G_H1=G/H$, the set of orbits of $G$ under
its right translation action by $H$. Hence, for any finite set
$\Sigma$ with continuous $H$-action, the unique map
$t_\Sigma:\Sigma\to 1$ yields a $G$-equivariant map
$$
\Ind^G_Ht_\Sigma:\Ind^G_H\to G/H
$$
and therefore $\Ind^G_H$ factors through a functor
\set\begin{equation}\label{eq_interposed}
H\tdu\Set\to G\tdu\Set/(G/H).
\end{equation}
It is easily seen that \eqref{eq_interposed} is an equivalence.
Indeed, one obtains a natural quasi-inverse, by the rule :
$(f:\Sigma\to G/H)\mapsto f^{-1}(H)$. The detailed verification
shall be left to the reader.

\begin{definition}\label{def_Galois-fibre-fct}
(\cite[Exp.V, Def.5.1]{SGA1}) Let $\cC$ be a category, and
$F:\cC\to\Set$ a functor.
\begin{enumerate}
\item
We say that $\cC$ is a {\em Galois category}, if $\cC$ is
equivalent to $P\tdu\Set$, for some profinite group $P$.
We denote $\Galois$ the category whose objects are all the
Galois categories, and whose morphisms are the exact functors
between Galois categories.
\item
We say that $F$ is a {\em fibre functor}, if $F$ is exact
and conservative, and $F(X)$ is a finite set for every
$X\in\Ob(\cC)$.
\item
We denote $\mathbf{fibre.Fun}$ the {\em $2$-category of fibre
functors}, defined as follows :
\begin{enumerate}
\item
The objects are all the pairs $(\cC,F)$ consisting of a
Galois category $\cC$ and a fibre functor $F$ for $\cC$.
\item
The $1$-cells $(\cC_1,F_1)\to(\cC_2,F_2)$ are all the pairs
$(G,\beta)$ consisting of an exact functor $G:\cC_1\to\cC_2$
and an isomorphism of functors $\beta:F_1\isom F_2\circ G$.
\item
And for every pair of $1$-cells
$(G',\beta'),(G,\beta):(\cC_1,F_1)\to(\cC_2,F_2)$, the
$2$-cells $(G',\beta')\to(G,\beta)$ are the isomorphisms
$\gamma:G'\isom G$ such that $(F_2*\gamma)\odot\beta'=\beta$.
\end{enumerate}
Composition of $1$-cells and $2$-cells is defined in the
obvious way. We shall also denote simply by $G$ a $1$-cell
$(G,\beta)$ as in (b), such that $F_1=F_2\circ G$ and $\beta$
is the identity automorphism of $F_1$.
\end{enumerate}
\end{definition}

\sset\subsubsection{}\label{subsec_Galois-cat-fib}
Notice that any Galois category $\cC$ admits a fibre functor :
indeed, if $P$ is any profinite group, the forgetful functor
$$
\sff_P:P\tdu\Set\to\Set
$$
fulfills the conditions of definition \ref{def_Galois-fibre-fct}(ii),
therefore the same holds for the functor $\sff_P\circ\beta$, if
$\beta:\cC\to P\tdu\Set$ is any equivalence. For any Galois
category $\cC$ and any fibre functor $F:\cC\to\Set$, we denote
$$
\pi_1(\cC,F)
$$
the group of automorphisms of $F$, and we call it the
{\em fundamental group of\/ $\cC$ pointed at $F$}. By definition,
for every $X\in\Ob(\cC)$, the finite set $F(X)$ is endowed
with a natural left action of $\pi_1(\cC,F)$. For every
$X\in\Ob(\cC)$ and every $\xi\in F(X)$, the stabilizer
$H_{X,\xi}\subset\pi_1(\cC,F)$ is a subgroup of finite index,
and we endow $\pi_1(\cC,F)$ with the coarsest group topology
for which all such $H_{X,\xi}$ are open subgroups.
The resulting topological group $\pi_1(\cC,F)$ is profinite,
and its natural left action on every $F(X)$ is continuous.
Thus, $F$ upgrades to a functor denoted
$$
F^\dagger:\cC\isom\pi_1(\cC,F)\tdu\Set.
$$
A basic result states that $F^\dagger$ is an equivalence
(\cite[Exp.V, Th.4.1]{SGA1}).

\begin{example}\label{ex_recover-profinite}
Let $P$ be any profinite group. Then there is an obvious
injective map
$$
P\to\pi_1(P\tdu\Set,\sff_P)
$$
and \cite[Exp.V, Th.4.1]{SGA1} implies that this map is an
isomorphism of profinite groups. In other words, the group
$P$ can be recovered, up to unique isomorphism, from the
category $P\tdu\Set$ together with its forgetful functor
$\sff_P$.
\end{example}

\begin{remark}\label{rem_funct-Galoisienne}
Let $\cC$, $\cC'$ be two Galois categories, and $F:\cC\to\Set$
a fibre functor.

(i)\ \
Any exact functor $G:\cC'\to\cC$ induces a continuous group
homomorphism :
$$
\pi_1(G):\pi_1(\cC,F)\to\pi_1(\cC',F\circ G)
\qquad
\omega\mapsto\omega*G.
$$

(ii)\ \
Furthermore, any isomorphism $\beta:F'\isom F$ of fibre functors
of $\cC$ induces an isomorphism of profinite groups :
$$
\pi_1(\beta):\pi_1(\cC',F)\isom\pi_1(\cC',F)
\qquad
\omega\mapsto\beta^{-1}\odot\omega\odot\beta
$$
(see \cite[Exp.V, \S4]{SGA1} for all these generalities).

(iii)\ \
Let now $(G,\beta):(\cC_1,F_1)\to(\cC_2,F_2)$ be a $1$-cell
of $\mathbf{fibre.Fun}$. Combining (i) and (ii), we deduce a
natural continuous group homomorphism
$$
\pi_1(G,\beta):\pi_1(\cC_2,F_2)\xrightarrow{\ \pi_1(G)\ }
\pi_1(\cC_1,F_2\circ G)\xrightarrow{\ \pi_1(\beta)\ }\pi_1(\cC_1,F_1).
$$
\end{remark}

\begin{proposition}\label{prop_nice-align}
With the notation of remark {\em\ref{rem_funct-Galoisienne}},
the rule that assigns :
\begin{itemize}
\item
To any object $(\cC,F)$ of\/ $\mathbf{fibre.Fun}$, the
profinite group $\pi_1(\cC,F)$
\item
To any $1$-cell $(G,\beta)$ of\/ $\mathbf{fibre.Fun}$, the
continuous map $\pi_1(G,\beta)$
\end{itemize}
defines a pseudo-functor
$$
\pi_1:\mathbf{fibre.Fun}\to\mathbf{pf.Grp}^o
$$
from the $2$-category of fibre functors, to the opposite
of the category of profinite groups (and continuous group
homomorphisms).
\end{proposition}
\begin{proof} (Here we regard $\mathbf{pf.Grp}$ as a
$2$-category with trivial $2$-cells : see example
\ref{ex_2-cats}(i)). Let
$$
(\cC_1,F_1)\xrightarrow{\ (G,\beta_G)\ }(\cC_2,F_2)
\xrightarrow{\ (H,\beta_H)\ }(\cC_3,F_3)
$$
be any pair of (composable) $1$-cells; functoriality on
$1$-cells amounts to the identity :
$$
\pi_1(G,\beta_G)\circ\pi_1(H,\beta_H)=
\pi_1(H\circ G,(\beta_H*G)\circ\beta_G)
$$
whose detailed verification we leave to the reader. Next,
let $\gamma:(G',\beta')\to(G,\beta)$ be a $2$-cell between
$1$-cells $(G',\beta'),(G,\beta):(\cC_1,F_1)\to(\cC_2,F_2)$;
we have to check that $\pi_1(G,\beta)=\pi_1(G',\beta')$.
This identity boils down to the commutativity of the diagram :
$$
\xymatrix{
\pi_1(\cC_2,F_2) \ar[rr]^-{\pi_1(G')} \ar[d]_{\pi_1(G)} & &
\pi_1(\cC_1,F_2\circ G') \ar[d]^{\pi_1(\beta')} \\
\pi_1(\cC_2,F_2\circ G) \ar[rr]^-{\pi_1(\beta)}
\ar[rru]^{\pi_1(F_2*\gamma)} & & \pi_1(\cC_1,F_1).
}$$
However, the commutativity of the lower triangular
subdiagram is clear, hence we are reduced to checking
the commutativity of the upper triangular subdiagram;
the latter is a special case of the following more
general :

\begin{claim} Let $\cC$ and $\cC'$ be two Galois categories,
$G,G':\cC'\to\cC$ two exact functors, $\beta:G'\isom G$ an
isomorphism, and $F:\cC\to\Set$ a fibre functor. Then the
induced diagram of profinite groups
$$
\xymatrix{
& \pi_1(\cC,F) \ar[ld]_{\pi_1(G)} \ar[rd]^{\pi_1(G')} \\
\pi_1(\cC',F\circ G) \ar[rr]^-{\pi_1(F*\beta)} & &
\pi_1(\cC',F\circ G')
}$$
commutes.
\end{claim}
\begin{pfclaim}[] Left to the reader.
\end{pfclaim}
\end{proof}

\begin{example}\label{ex_recover-mapfrom-ind}
Let $\omega:P\to Q$ be a continuous group homomorphism between
profinite groups. Clearly $\sff_P\circ\Res(\omega)=\sff_Q$, and
it is easily seen that the resulting diagram
$$
\xymatrix{ P \ar[rr]^-\omega \ar[d] & & Q \ar[d] \\
\pi_1(P\tdu\Set,\sff_P) \ar[rr]^-{\pi_1(\Res(\omega))}
& & \pi_1(Q\tdu\Set,\sff_Q)
}$$
commutes, where the vertical arrows are the natural identifications
given by example \ref{ex_recover-profinite} : the verification
is left as an exercise to the reader.
\end{example}

\sset\subsubsection{}\label{subsec_profin-lims}
Let $\underline P:=(P_i~|~i\in I)$ be a cofiltered system of
profinite groups, with continuous transition maps, and denote
by $P$ the limit of this system, in the category of groups.
Then $P$ is naturally a closed subgroup of $Q:=\prod_{i\in I}P_i$,
and the topology $\cT$ induced by the inclusion map $P\to Q$
makes it into a compact and complete topological group. Moreover,
since the topology of $Q$ is profinite, the same holds for the
topology $\cT$ (details left to the reader). It is then easily
seen that the resulting topological group $(P,\cT)$ is the limit
of the system $\underline P$ in the category of profinite groups.

\begin{proposition}\label{prop_look-a-prop}
In the situation of \eqref{subsec_profin-lims}, the natural
functor
$$
\Pscolim{i\in I}P_i\tdu\Set\to P\tdu\Set
$$
is an equivalence.
\end{proposition}
\begin{proof} The functor is obviously faithful; let us show
that it is also full. Indeed, let $j\in I$ be any index,
$\Sigma,\Sigma'$ two objects of $P_j\tdu\Set$, and
$\phi:\Sigma\to\Sigma'$ a $P$-equivariant map; we need
to show that $\phi$ is already $P_i$-equivariant, for
some index $i\in I$. To this aim, we may as usual replace
$I$ by $I/j$, and assume that $j$ is the final element of
$I$. We may also find an open normal subgroup $H_j\subset P_j$
that acts trivially on both $\Sigma$ and $\Sigma'$. Then,
for every $i\in I$, we let $H_i\subset P_i$ be the preimage
of $P_j$, and we set $\bar P_i:=P_i/H_i$. Let also
$\bar P:=P/H$, where $H\subset P$ is the preimage of $H_j$;
by construction, $\phi$ is $\bar P$-equivariant. Clearly,
we may find $i\in I$ such that the image of $\bar P$ in
the finite group $\bar P_j$ equals the image of $\bar P_i$,
and for such index $i$, the induced map $\bar P\to\bar P_i$
is an isomorphism. Especially, $\phi$ is $P_i$-equivariant,
as sought.

Lastly, we show essential surjectivity. Indeed, let $\Sigma$
be any object of $P\tdu\Set$; we have to show that the
$P$-action on $\Sigma$ is the restriction of a continuous
$P_i$-action, for a suitable $i\in I$. However, we may find
a normal open subgroup $H\subset P$ that acts trivially on
$\Sigma$. Then there exists a normal open subgroup $L\subset Q$
(notation of \eqref{subsec_profin-lims}) such that
$P\cap L\subset H$. We may also assume that there exist
a finite subset $J\subset I$ and for every $i\in J$ an
open normal subgroup $L_i\subset P_i$ such that
$L=\prod_{i\in J}L_i\times\prod_{i\in I\setminus J}P_i$.
Pick an index $j\in I$ that admits morphisms $f_i:j\to i$
in $I$, for every $i\in J$, and let $L'_i\subset P_j$
denote the preimage of $L_i$ under the corresponding
map $P_j\to P_i$. Finally, set $H_j:=\bigcap_{i\in J}L'_i$.
By construction, $H$ contains the preimage of $H_j$ in $P$,
and we may therefore assume that $H$ is this preimage.
We may replace as usual $I$ by $I/j$, and assume that $j$
is the final element of $I$. Then, for every $i\in I$, we
let $H_i$ denote the preimage of $H_j$ in $P_i$, and we
set $\bar P_i:=P_i/H_i$. Set as well $\bar P:=P/H$.
Clearly, there exists $i\in I$ such that the image of
the induced map $\bar P\to\bar P_j$ equals the image
of $\bar P_i$; for such index $i$, the induced map
$\bar P\to\bar P_i$ is an isomorphism. Thus, $\Sigma$
the restriction of an object of $P_i\tdu\Set$, as wished.
\end{proof}

\sset\subsubsection{}\label{subsec_profin-pseudolim}
We consider now a situation that generalizes slightly that
of \eqref{subsec_profin-lims}. Namely, let $I$ be a small
filtered category, and
$$
(\cC_\bullet,F_\bullet):I\to\mathbf{fibre.Fun}
\qquad
i\mapsto(\cC_i,F_i)
$$
a pseudo-functor. By proposition \ref{prop_nice-align},
the composition of $\pi_1$ and $(\cC_\bullet,F_\bullet)$
is a functor
$$
\pi_1(\cC_\bullet,F_\bullet):I^o\to\mathbf{pf.Grp}
\qquad
i\mapsto P_i:=\pi_1(\cC_i,F_i).
$$
Let $P$ denote the limit (in $\mathbf{pf.Grp}$) of the
cofiltered system $P_\bullet$, and set
$$
\cC:=\Pscolim{I}\cC_\bullet
$$
where the $2$-colimit is formed in the $2$-category of
small categories. We may then state :

\begin{corollary}\label{cor_profin-pseudolim}
In the situation of \eqref{subsec_profin-pseudolim}, there
exists a natural equivalence :
$$
\cC\isom P\tdu\Set.
$$
\end{corollary}
\begin{proof} Recall that $(\cC_\bullet,F_\bullet)$ is the
datum of isomorphisms
$$
\beta_\phi:F_j\isom F_i\circ\cC_\phi
\qquad
\text{for every morphism $\phi:j\to i$ in $I$}
$$
and $2$-cells :
$$
\tau_{\psi,\phi}:(\cC_{\psi\circ\phi},\beta_{\psi\circ\phi})
\isom(\cC_\psi,\beta_\psi)\circ(\cC_\phi,\beta_\phi)
\qquad
\text{for every composition
$j\xrightarrow{\ \phi\ }i\xrightarrow{\ \psi\ }k$}
$$
that -- by definition -- satisfy the identities :
\set\begin{equation}\label{eq_not-coherence}
(\beta_\psi*\cC_\phi)\odot\beta_\phi=
(F_k*\tau_{\psi,\phi})\odot\beta_{\psi\circ\phi}
\qquad
\text{for every composition
$j\xrightarrow{\ \phi\ }i\xrightarrow{\ \psi\ }k$}
\end{equation}
as well as the composition identities :
$$
((\cC_\mu,\beta_\mu)*\tau_{\psi,\phi})
\odot\tau_{\mu,\psi\circ\phi}=
(\tau_{\mu,\psi}*(\cC_\phi,\beta_\phi))
\odot\tau_{\mu\circ\psi,\phi}
\qquad
\text{for compositions $j\xrightarrow{\ \phi\ }
i\xrightarrow{\ \psi\ }k \xrightarrow{\ \mu\ }l$}.
$$
Let $P_\bullet\tdu\Set:I\to\bCat$ denote
the functor given by the rule : $i\mapsto P_i\tdu\Set$
for every $i\in\Ob(I)$, and $\phi\mapsto\Res(P_\phi)$, where
$P_\phi:=\pi_1(\cC_\phi,\beta_\phi)$ for every morphism
$\phi$ of $I$. In view of proposition \ref{prop_look-a-prop}
and lemma \ref{lem_pseudo-trivial}, it suffices to show that
the rule : $i\mapsto F^\dagger_i$ for every $i\in\Ob(I)$
(notation of \eqref{subsec_Galois-cat-fib}), extends to a
pseudo-natural isomorphism
$$
F^\dagger_\bullet:\cC_\bullet\isom P_\bullet\tdu\Set.
$$
Indeed, let $\phi:j\to i$ be any morphism of $I$, and $X$
any object of $\cC_j$; we remark that the bijection
$$
\beta_\phi(X):
\Res(P_\phi)(F^\dagger_jX)\isom F^\dagger_i\circ\cC_\phi(X)
$$
is $P_i$-equivariant; the proof amounts to unwinding the
definitions, and shall be left to the reader. Hence we
get an isomorphism of functors
$\beta^\dagger_\phi:
\Res(P_\phi)\circ F^\dagger_j\isom F_i^\dagger\circ\cC_\phi$,
and from \eqref{eq_not-coherence} we deduce the identities :
$$
(\beta^\dagger_\psi*\cC_\phi)\odot
(\Res(P_\phi)*\beta^\dagger_\phi)=(F^\dagger_k*\tau_{\psi,\phi})
\odot\beta^\dagger_{\psi\circ\phi}
\qquad
\text{for every composition
$j\xrightarrow{\ \phi\ }i\xrightarrow{\ \psi\ }k$}.
$$
The latter show that the system $\beta^\dagger_\bullet$ fulfills
the coherence axiom for a pseudo-natural transformation, as
required.
\end{proof}

\begin{remark}\label{rem_2-colim-Galois}
(i)\ \
Keep the situation of \eqref{subsec_profin-pseudolim}, and
let $a_\bullet:\cC_\bullet\Rightarrow\cC$ be the universal
pseudo-cocone induced by the pseudo-functor $\cC_\bullet$.
We may regard the pseudo-functor $(\cC_\bullet,F_\bullet)$
as a pseudo-cocone $F_\bullet:\cC_\bullet\Rightarrow\Set$
whose vertex is the category $\Set$. Then, by the universal
property of colimits, we get a functor $F:\cC\to\Set$ and
an isomorphism
$$
\sigma_\bullet:F*a_\bullet\isom F_\bullet.
$$
On the other hand, let
$r_\bullet:P_\bullet\tdu\Set\Rightarrow P\tdu\Set$
be the natural cocone (so $r_i$ is the restriction functor
corresponding to the natural map $P\to P_i$, for every
$i\in\Ob(I)$); the equivalence $G$ of corollary
\ref{cor_profin-pseudolim} is deduced from the pseudo-cocone
$r_\bullet\circ F^\dagger_\bullet:\cC_\bullet\to P\tdu\Set$
(where $F^\dagger_\bullet$ is as in the proof of corollary
\ref{cor_profin-pseudolim}), so we have an isomorphism of
pseudo-functors
$$
t_\bullet:G*a_\bullet\isom r_\bullet\circ F^\dagger_\bullet
$$
whence an isomorphism
$$
\sff_P*t_\bullet:(\sff_P\circ G)*a_\bullet\isom
\sff_P*(r_\bullet\circ F^\dagger_\bullet)=F_\bullet
$$
(notation of \eqref{subsec_Galois-cat-fib}). There follows an
isomorphism $F*a_\bullet\isom(\sff_P\circ G)*a_\bullet$;
by the universal property of the $2$-colimit, the latter must
come from a unique isomorphism $\theta:F\isom\sff_P\circ G$.
Especially, we see that $F$ is also a fibre functor, and we
get a pseudo-cocone
$$
(a_\bullet,\sigma_\bullet):
(\cC_\bullet,F_\bullet)\Rightarrow(\cC,F).
$$
It is now immediate that $(\cC,F)$ is the $2$-colimit (in
the $2$-category $\mathbf{fibre.Fun}$) of the pseudo-functor
$(\cC_\bullet,F_\bullet)$, and $(a_\bullet,\sigma_\bullet)$
is the corresponding universal pseudo-cocone.

(ii)\ \
Likewise, $r_\bullet$ may be regarded as a universal
pseudo-cocone
$$
r_\bullet:(P_\bullet\tdu\Set,\sff_{P_\bullet})
\Rightarrow(P\tdu\Set,\sff_P)
$$
(with trivial coherence constraint), and the coherence
constraint $\beta^\dagger_\bullet$ as in the proof of
corollary \ref{cor_profin-pseudolim} yields a pseudo-natural
equivalence
$$
F^\dagger_\bullet:(\cC_\bullet,F_\bullet)\isom
(P_\bullet\tdu\Set,\sff_{P_\bullet})
$$
as well as an isomorphism
$$
(G,\theta)*(a_\bullet,\sigma_\bullet)\isom
r_\bullet\odot F^\dagger_\bullet.
$$
Thus, for every $i\in\Ob(I)$ we get a $2$-cell of
$\mathbf{fibre.Fun}$ :
$$
(G,\theta)\circ(a_i,\sigma_i)\to r_i\circ F^\dagger_i
$$
whence -- by proposition \ref{prop_nice-align} -- a commutative
diagram of profinite groups :
$$
\xymatrix{ P \ar[rr] \ar[d]_{\pi_1(G,\theta)} & &
P_i \ar[d]^{\pi_1(F^\dagger_i)} \\
\pi_1(\cC,F) \ar[rr]^-{\pi_1(a_i,\sigma_i)} & &
\pi_1(\cC_i,F_i).
}$$
\end{remark}

\sset\subsubsection{}\label{subsec_sub-Galois}
Let $(\cC,F)$ be a fibre functor, and $X$ a connected
object of $\cC$ (example \ref{ex_equalizers}(iv)); pick
any $\xi\in F(X)$, and let $H_\xi\subset\pi_1(\cC,F)$ be
the stabilizer of $\xi$ for the natural left action of
$\pi_1(\cC,F)$ on $F(X)$. For every object $f:Y\to X$ of
$\cC\!/\!X$, we set
$$
F_\xi(f):=F(f)^{-1}(\xi)\subset F(Y).
$$
It is clear that the rule $f\mapsto F_\xi(f)$ yields a functor
$F^\dagger_\xi:\cC\!/\!X\to H_\xi\tdu\Set$, which we call
the {\em subfunctor of $F_{|X}$ selected by $\xi$}.

\begin{proposition}\label{prop_sarah-call}
In the situation of \eqref{subsec_sub-Galois}, we have :
\begin{enumerate}
\item
$\cC\!/\!X$ is also a Galois category, and
$F_\xi:=\sff_{H_\xi}\circ F^\dagger_\xi$ is a fibre functor
for $\cC\!/\!X$.
\item
The functor $F^\dagger_\xi$ induces a natural isomorphism of
profinite groups :
$$
\pi_1(F^\dagger_\xi):H_\xi\isom\pi_1(\cC\!/\!X,F_\xi).
$$
\end{enumerate}
\end{proposition}
\begin{proof} The fibre functor $F$ induces an equivalence
of categories
$$
\cC\!/\!X\isom\pi_1(\cC,F)\tdu\Set/F(X).
$$
On the other hand, since $X$ is connected, there exists a
unique isomorphism $\omega:F(X)\isom G/H_\xi$ of $G$-sets
such that $\omega(\xi)=H_\xi$, and then the discussion of
\eqref{subsec_nat-quinverse} yields an equivalence
$$
\pi_1(\cC,F)\tdu\Set/F(X)\isom H_\xi\tdu\Set.
$$
A simple inspection shows that the resulting equivalence
$\cC\!/\!X\isom H_\xi\tdu\Set$ is none else than the
functor $F^\dagger_\xi$, so the assertion follows from
remark \ref{rem_funct-Galoisienne}(i) and example
\ref{ex_recover-profinite}.
\end{proof}

\sset\subsubsection{}\label{subsec_cleavage-motoko}
Let $(\cC,F)$ be a fibre functor, and let us now fix a
cleavage $\sc:\cC^o\to\bCat$ for the fibred category
$\st:\sMorph(\cC)\to\cC$ (see example
\ref{ex_fibred-cats}(iii)). Also, let $I$ be a small
cofiltered category, and $X_\bullet:I\to\cC$ a functor
such that $X_i$ is a connected object of $\cC$, for
every $i\in\Ob(I)$; we pick an element
$$
\xi_\bullet\in\lim_I F\circ X_\bullet.
$$
In other words, $\xi_\bullet:=(\xi_i\in F(X_i)~|~i\in I)$
is a compatible system of elements such that
$$
F(\phi)(\xi_j)=\xi_i
\qquad
\text{for every morphism $\phi:j\to i$ in $I$.}
$$
For every $i\in I$, we denote by $H_i\subset\pi_1(\cC,F)$ the
stabilizer of $\xi_i$ for the left action of $\pi_1(\cC,F)$
on $F(X_i)$. Clearly, any morphism $j\to i$ induces an
inclusion $H_j\subset H_i$. Furthermore, let
$$
F^\dagger_i:\cC\!/\!X_i\to H_i\tdu\Set
\qquad
\text{for every $i\in I$}
$$
be the subfunctor selected by $\xi_i$, and set
$F_i:=\sff_{H_i}\circ F^\dagger_i$. Let $\phi:j\to i$
be a morphism of $I$; to the corresponding morphism
$X_\phi:X_j\to X_i$, the cleavage $\sc$ attaches a
pull-back functor
$$
X^*_\phi:\cC\!/\!X_i\to\cC\!/\!X_j.
$$
Especially, for any object $Y\in\Ob(\cC/X_i)$ we have the
cartesian diagram in $\cC$ :
$$
\xymatrix{ X^*_\phi(Y) \ar[r] \ar[d] & Y \ar[d] \\
           X_j \ar[r]^-{X_\phi} & X_i
}$$
whence, since $F$ is exact, a natural bijection :
$$
F(X^*_\phi(Y))\isom F(Y)\times_{F(X_i)}F(X_j)
$$
which in turns yields a bijection :
$$
F^\dagger_j(X^*_\phi(Y))\isom F(Y)\times_{F(X_i)}\{\xi_j\}=
F^\dagger_i(Y)\times\{\xi_j\}.
$$
That is, we have a natural isomorphism of functors :
\set\begin{equation}\label{eq_restrict-action}
\alpha_\phi^\dagger:F^\dagger_j\circ X^*_\phi\isom
\Res^{H_j}_{H_i}\circ F^\dagger_i
\end{equation}
and since the $X^*_\phi$ are exact functors, it is easily
seen that the isomorphisms
$\alpha_\phi:=\sff_j*\alpha^\dagger_\phi$ yield a pseudo-functor
\set\begin{equation}\label{eq_labeyrie}
I\to\mathbf{fibre.Fun}
\qquad
i\mapsto(\cC\!/\!X_i,F_i)
\qquad
\phi\mapsto(X^*_\phi,\alpha_\phi).
\end{equation}
Moreover, $\alpha_\phi^\dagger$ can be seen as a $2$-cell of
$\mathbf{fibre.Fun}$ : $F^\dagger_j\circ(X^*_\phi,\alpha_\phi)
\isom\Res^{H_i}_{H_j}\circ F^\dagger_i$, whence a commutative
diagram of profinite groups :
$$
\xymatrix{ H_j \ar[rr] \ar[d]_{\pi_1(F^\dagger_j)} & &
H_i \ar[d]^{\pi_1(F^\dagger_i)} \\
\pi_1(\cC\!/\!X_j,F_j) \ar[rr]^-{\pi_1(X^*_\phi,\alpha_\phi)}
& & \pi_1(\cC\!/\!X_i,F_i)
}$$
whose top horizontal arrow is the inclusion map. Especially,
notice that the map $\pi_1(X^*_\phi,\alpha_\phi)$ does not
depend on the chosen cleavage; this can also be seen by
remarking that any two cleavages $\sc$, $\sc'$ are related
by a pseudo-natural isomorphism $\sc\isom\sc'$ (details left
to the reader).

\sset\subsubsection{}\label{subsec_intersect-Galois}
Let $(\cC\!/\!X,F_\xi)$ be the $2$-colimit of the pseudo-functor
\eqref{eq_labeyrie}, as in remark \ref{rem_2-colim-Galois}(i),
and fix a corresponding universal pseudo-cocone
$(a_\bullet,\sigma_\bullet):(\cC\!/\!X_\bullet,F_\bullet)\Rightarrow
(\cC\!/\!X,F_\xi)$. We may then state :

\begin{corollary}\label{cor_intersect-Hs}
In the situation of \eqref{subsec_intersect-Galois},
there exists a natural isomorphism of profinite groups :
$$
H:=\bigcap_{i\in\Ob(I)}H_i\isom\pi_1(\cC\!/\!X,F_\xi)
$$
which fits into a commutative diagram :
\set\begin{equation}\label{eq_no-espace}
{\diagram H \ar[rr] \ar[d] & & H_i \ar[d]^{\pi_1(F^\dagger_i)} \\
\pi_1(\cC\!/\!X,F_\xi) \ar[rr]^-{\pi_1(a_i,\sigma_i)} & &
\pi_1(\cC\!/\!X_i,F_i)
\enddiagram}
\qquad
\text{for every $i\in\Ob(I)$}
\end{equation}
whose top horizontal arrow is the inclusion map.
\end{corollary}
\begin{proof} By corollary \ref{cor_profin-pseudolim}, we have
an isomorphism of $\pi_1(\cC\!/\!X,F_\xi)$ with the limit of
the cofiltered system $(\pi_1(\cC\!/\!X_,F_i)~|~i\in\Ob(I))$;
on the other hand, the discussion of \eqref{subsec_cleavage-motoko}
shows that the latter system is naturally isomorphic to the
system $(H_i~|~i\in\Ob(I))$. Lastly, the commutativity of
\eqref{eq_no-espace} follows from remark
\ref{rem_2-colim-Galois}(ii).
\end{proof}

\subsection{Tensor categories and abelian categories}
\label{sec_tensors-ab}
In this section we assemble some basic definitions and
results that pertain to abelian categories and other
related classes of categories with extra structure.

\begin{definition}\label{def_tensor-cat}
A {\em tensor category} is a datum
$\underline\cC:=(\cC,\otimes,\Phi,\Psi)$ consisting of
a category $\cC$, a functor
$$
\otimes:\cC\times\cC\to\cC
\quad :\quad
(X,Y)\mapsto X\otimes Y
$$
and natural isomorphisms :
$$
\Phi_{X,Y,Z}:X\otimes(Y\otimes Z)\isom(X\otimes Y)\otimes Z
\qquad
\Psi_{X,Y}:X\otimes Y\isom Y\otimes X
$$
for every $X,Y,Z\in\Ob(\cC)$, called respectively the
{\em associativity} and {\em commutativity constraints}
of $\underline\cC$, that satisfy the following axioms.
\begin{enumerate}
\alphaenu
\item
{\em Coherence axiom} : the diagram
$$
\xymatrix{ X\otimes(Y\otimes(Z\otimes T))
\ar[rr]^-{\Phi_{X,Y,Z\otimes T}} \ar[d]_{X\otimes\Phi_{Y,Z,T}} & &
(X\otimes Y)\otimes(Z\otimes T) \ar[rr]^-{\Phi_{X\otimes Y,Z,T}} & &
((X\otimes Y)\otimes Z)\otimes T \\
X\otimes((Y\otimes Z)\otimes T)
\ar[rrrr]^-{\Phi_{X,Y\otimes Z,T}} & & & &
(X\otimes(Y\otimes Z))\otimes T \ar[u]_{\Phi_{X,Y,Z}\otimes T}
}$$
commutes, for every $X,Y,Z,T\in\Ob(\cC)$.
\item
{\em Compatibility axiom} : the diagram
$$
\xymatrix{ X\otimes(Y\otimes Z) \ar[rr]^-{\Phi_{X,Y,Z}}
\ar[d]_{X\otimes\Psi_{Y,Z}}
& & (X\otimes Y)\otimes Z \ar[rr]^-{\Psi_{X\otimes Y,Z}} & &
Z\otimes(X\otimes Y) \ar[d]^{\Phi_{Z,X,Y}} \\
X\otimes(Z\otimes Y) \ar[rr]^-{\Phi_{X,Z,Y}} & &
(X\otimes Z)\otimes Y \ar[rr]^-{\Psi_{X,Z}\otimes Y} & &
(Z\otimes X)\otimes Y
}$$
commutes, for every $X,Y,Z\in\Ob(\cC)$.
\item
{\em Commutation axiom} : we have
$\Psi_{Y,X}\circ\Psi_{X,Y}=\one_{X\otimes Y}$ for every
$X,Y\in\Ob(\cC)$
\item
{\em Unit axiom} : there exist an object $U\in\Ob(\cC)$
and an isomorphism $u:U\isom U\otimes U$ such that the functor
\set\begin{equation}\label{eq_auto-morphism}
\cC\to\cC
\quad :\quad
X\mapsto U\otimes X
\end{equation}
is an equivalence. One says that $(U,u)$ is a
{\em unit object} of $\underline\cC$.
\end{enumerate}
\end{definition}

\begin{lemma}\label{lem_second-compat}
Let $\underline\cC:=(\cC,\otimes,\Phi,\Psi)$ be any tensor
category. The diagram
$$
\xymatrix{ X\otimes(Y\otimes Z) \ar[rr]^-{\Phi_{X,Y,z}}
\ar[d]_{X\otimes\Psi_{Y,Z}} & &
(X\otimes Y)\otimes Z \ar[rr]^-{\Psi_{X,Y}\otimes Z} & &
(Y\otimes X)\otimes Z \\
X\otimes(Z\otimes Y) \ar[rr]^-{\Phi_{X,Z,Y}} & &
(X\otimes Z)\otimes Y \ar[rr]^-{\Psi_{X\otimes Z,Y}} & &
Y\otimes(X\otimes Z) \ar[u]_{\Phi_{Y,X,Z}}
}$$
commutes, for every $X,Y,Z\in\Ob(\cC)$.
\end{lemma}
\begin{proof} To ease notation, we shall omit the tensor
symbol $\otimes$ between objects, and we shall drop the
subscript from $\Phi$ and $\Psi$ when we display a
diagram. It suffices to consider the diagram
$$
\xymatrix{
& Y(XZ) \ar[rrrr]^-\Phi \ar[d]^{Y\otimes\Psi} & & &
& (YX)Z \ar[d]_{\Psi} \\
& Y(ZX) \ar[r]^-\Phi &  (YZ)X \ar[rr]^-{\Psi\otimes X}
& & (ZY)X  & Z(YX) \ar[l]_\Phi \\
(XZ)Y \ar[r]^{\Psi\otimes Y} \ar[ruu]^\Psi
& (ZX)Y \ar[u]_\Psi & & & & Z(XY) \ar[u]^{Z\otimes\Psi}
\ar[llll]_{\Phi} & (XY)Z \ar[l]_{\Psi} \ar[luu]_{\Psi\otimes Z} \\
X(ZY) \ar[u]^\Phi & & & & & &
\ar[llllll]_-{X\otimes\Psi} X(YZ) \ar[u]_\Phi
}$$
whose two triangular subdiagrams commute by naturality of
$\Psi$, and whose three rectangular subdiagrams commute by
compatibility.
\end{proof}

\begin{definition} Let $\underline\cC:=(\cC,\otimes,\Phi,\Psi)$
and $\underline\cC':=(\cC',\otimes',\Phi',\Psi')$ be two tensor
categories. A {\em tensor functor} $\underline\cC\to\underline\cC'$
is a pair $(F,c)$ consisting of a functor $F:\cC\to\cC'$ and
a natural isomorphism
$$
c_{X,Y}:FX\otimes FY\isom F(X\otimes Y)
\qquad
\text{for all $X,Y\in\Ob(\cC)$}
$$
such that the following holds.
\begin{enumerate}
\alphaenu
\item
For every objects $X,Y,Z$ of $\cC$, the diagram
$$
\xymatrix{ FX\otimes(FY\otimes FZ)
\ar[rr]^-{FX\otimes c_{YZ}} \ar[d]_{\Phi'_{FX,FY,FZ}} & &
FX\otimes F(Y\otimes Z) \ar[rr]^-{c_{X,Y\otimes Z}} & &
F(X\otimes(Y\otimes Z)) \ar[d]^{F(\Phi_{X,Y,Z})} \\
(FX\otimes FY)\otimes FZ \ar[rr]^-{c_{X,Y}\otimes FZ} & &
F(X\otimes Y)\otimes FZ \ar[rr]^-{c_{X\otimes Y,Z}} & &
F((X\otimes Y)\otimes Z)
}$$
commutes.
\item
For all objects $X$, $Y$ of $\cC$, the diagram
$$
\xymatrix{ FX\otimes FY \ar[rr]^-{c_{X,Y}} \ar[d]_{\Psi_{FX,FY}} & &
F(X\otimes Y) \ar[d]^{F(\Psi_{X,Y})} \\
FY\otimes FX \ar[rr]^-{c_{Y,X}} & & F(Y\otimes X)
}$$
commutes.
\item
If $(U,u)$ is a unit object of $\underline\cC$, then
$(FU,c_{U,U}^{-1}\circ Fu)$ is a unit object of $\underline\cC'$.
\end{enumerate}
\end{definition}

\begin{remark}\label{rem_Mac-Lane}
(i)\ \
Lemma \ref{lem_second-compat} illustrates a general principle
valid in every tensor categor $\underline\cC$: namely, say
that $X_1,\dots,X_n$ is a sequence of {\em distinct} objects of
$\cC$, and $X'$ and $X''$ are obtained from these two sequences
by taking tensor products several times, and in any order, in
which case we say that $X'$ and $X''$ {\em have no repetitions}.
Now, there will be usually various ways to combine the
associativity and commutativity constraints, in order to
exhibit some isomorphism $X'\isom X''$. However, the
resulting isomorphism shall be independent of the way in
which it is expressed as such a combination. This follows
from a theorem of Mac Lane. To formalize this result, one
could observe that, for any set $\Sigma$, there exists a
universal tensor category $T_\Sigma$ ``generated by $\Sigma$'',
{\em i.e.} such that -- for any other tensor category
$\underline\cC$ -- any mapping $\Sigma\to\Ob(\cC)$ extends
uniquely, up to isomorphism, to a tensor functor
$T_\Sigma\to\cC$. Then Mac Lane's theorem says that, for
every set $\Sigma$, and every object $X\in\Ob(T_\Sigma)$
that has no repetitions, the group $\Aut_{T_\Sigma}(X)$ is
trivial. Instead of relying on such a general result, we
shall make {\em ad hoc} verifications, as in the proof of
lemma \ref{lem_second-compat} and of the forthcoming proposition
\ref{prop_unit-isom}. However, in view of this principle,
in the following we shall often omit a detailed description
of the isomorphism that we choose to connect two given
objects that are thus related : the reader will be able
in any case to produce one isomorphism, and the principle
says that these choices cannot be source of ambiguities.

(ii)\ \
Let $\cD$ be any category, and $\underline\cC$ a tensor
category. Then notice that $\bFun(\cD,\cC)$ inherits
from $\underline\cC$ a natural tensor category structure :
we leave to the reader the task of spelling out the
details. Moreover notice that if
$(F,c):\underline\cC{}_1\to\underline\cC{}_2$ and
$(F',c'):\underline\cC{}_2\to\underline\cC{}_3$ are
any two tensor functors between tensor categories, then
the composition
$$
(F'\circ F,(F'*c)\odot(c'*(F\times F))):
\underline\cC{}_1\to\underline\cC{}_3
$$
is again a tensor functor.
\end{remark}

\begin{proposition}\label{prop_unit-isom}
Let $(U,u)$ be a unit object of a tensor category
$\underline\cC$. Then there exists a unique natural
isomorphism
$$
u_X:X\isom U\otimes X
\qquad
\text{for every $X\in\Ob(\cC)$}
$$
such that $u_U=u$, and such that the diagrams
$$
\xymatrix{ X\otimes Y \ar[rr]^-{u_{X\otimes Y}}
\ar[drr]_{u_X\otimes Y} & &
U\otimes(X\otimes Y) \ar[d]^{\Phi_{U,X,Y}} &
X\otimes Y \ar[rr]^-{u_X\otimes Y} \ar[d]_{X\otimes u_Y} & &
(U\otimes X)\otimes Y \ar[d]^{\Psi_{U,X}\otimes Y} \\
& & (U\otimes X)\otimes Y &
X\otimes(U\otimes Y) \ar[rr]^-{\Phi_{X,U,Y}} & &
(X\otimes U)\otimes Y
}$$
commute for every $X,Y\in\Ob(\cC)$.
\end{proposition}
\begin{proof} Since \eqref{eq_auto-morphism} is an equivalence,
there exists a unique isomorphism $u_X$ fitting into the
commutative diagram
$$
\xymatrix{ UX \ar[rr]^-{u\otimes X} \ar[rrd]_{U\otimes u_X} & &
(UU)X \\
& & U(UX) \ar[u]_\Phi.
}$$
With this definition, the naturality of the rule : $X\mapsto u_X$
is clear. In order to check the commutativity of the first diagram,
it suffices to show that
\set\begin{equation}\label{eq_with-Theta}
(U\otimes\Phi_{U,X,Y})\circ(U\otimes u_{XY})=U\otimes(u_X\otimes Y).
\end{equation}
However, set $\Theta:=(\Phi_{U,U,X}\otimes Y)\circ\Phi_{U,UX,Y}$;
we have :
$$
\begin{aligned}
\Theta\circ(U\otimes(u_X\otimes Y))
& =(\Phi_{U,U,X}\otimes Y)\circ((U\otimes u_X)\otimes Y)\circ
\Phi_{U,X,Y} \\
& =((u\otimes X)\otimes Y)\circ\Phi_{U,X,Y} \\
& =\Phi_{UU,X,Y}\circ(u\otimes(XY)) \\
& =\Phi_{UU,X,Y}\circ\Phi_{U,U,XY}\circ(U\otimes u_{XY}) \\
& =\Theta\circ(U\otimes\Phi_{U,X,Y})\circ(U\otimes u_{XY})
\end{aligned}
$$
where the first and third identities hold by naturality of $\Phi$,
the second and fourth by the definition of $u_X$ and respectively
$u_{X\otimes Y}$, and the fifth by coherence. Since $\Theta$ is an
isomorphism, we get \eqref{eq_with-Theta}.
Next, in light of the foregoing, the second diagram
commutes if and only if the diagram
\set\begin{equation}\label{eq_monster}
{\diagram XY \ar[rr]^-{u_{XY}}
\ar[rrd]_{X\otimes u_Y} & & U(XY)
\ar[rr]^-\Phi & &
(UX)Y \ar[d]^{\Psi\otimes Y} \\
& & X(UY) \ar[rr]^-\Phi & &
(XU)Y
\enddiagram}
\end{equation}
commutes, and it suffices to check that $U\otimes\eqref{eq_monster}$
commutes. To this aim, we consider the diagram
$$
\xymatrix{ U(XY) \ar[rr]^-{U\otimes(X\otimes u_Y)}
\ar[dd]_{u\otimes(XY)} \ar[rd]^{U\otimes u_{XY}} & &
U(X(UY)) \ar[rr]^-{U\otimes\Phi} & &
U((XU)Y)
\ar[d]^{U\otimes(\Psi\otimes Y)} \\
& U(U(XY)) \ar[rrr]^-{U\otimes\Phi} \ar[ld]_\Phi & & &
U((UX)Y) \ar[d]^\Phi \\
(UU)(XY) \ar[rr]^-\Phi & &
((UU)X)Y & & (U(UX))Y \ar[ll]_{\Phi\otimes Y}.
}$$
whose lower subdiagram commutes by the coherence axiom
for $(U,U,X,Y)$, whose upper subdiagram is equivalent
to \eqref{eq_monster}, since $\Psi_{U,X}^{-1}=\Psi_{X,U}$,
and whose triangular subdiagram commutes by definition
of $u_{XY}$. Hence, we are reduced to showing that the
outer rectangular subdiagram of the above diagram commutes.
However, we have a commutative diagram :
$$
\xymatrix{
U(X(UY)) \ar[d]_\Phi & &
U(XY)
\ar[ll]_-{U\otimes(X\otimes u_Y)} \ar[rr]^-{u\otimes(XY)}
\ar[d]_\Phi & &
(UU)(XY) \ar[d]^\Phi \\
(UX)(UY) \ar[d]_{\Psi\otimes(UY)}
& & (UX)Y 
\ar[ll]_-{(UX)\otimes u_Y} \ar[rr]^-{(u\otimes X)\otimes Y} 
\ar[d]_-{\Psi\otimes Y} & &
((UU)X)Y \ar[d]^{\Psi\otimes Y} \\
(XU)(UY) & & (XU)Y
\ar[rr]^-{(X\otimes u)\otimes Y}
\ar[ll]_-{(XU)\otimes u_Y} & & (X(UU))Y \\
& & X(UY) \ar[lld]_-{X\otimes(U\otimes u_Y)}
\ar[rrd]^-{X\otimes(u\otimes Y)} \ar[u]^\Phi \\
X(U(UY)) \ar[rrrr]^-{X\otimes\Phi}
\ar[uu]^\Phi & & & &
X((UU)Y) \ar[uu]_\Phi
}$$
so we are further reduced to checking the commutativity of the
diagram :
$$
\xymatrix@C-1pt{
& U(X(UY)) \ar[r]^-{U\otimes\Phi}
\ar[ld]_\Phi &
U((XU)Y)
\ar[rr]^-{U\otimes(\Psi\otimes Y)} & &
U((UX)Y) \ar[rd]^\Phi \\
(UX)(UY) \ar[d]_{\Psi\otimes(UY)}
& U((UY)X) \ar[u]^{U\otimes\Psi} \ar[d]_\Phi
& U(U(YX))
\ar[l]_-{U\otimes\Phi} \ar[rr]^-{U\otimes(U\otimes\Psi)} \ar[d]_\Phi
& & U(U(XY)) \ar[u]^{U\otimes\Phi} \ar[d]_\Phi &
(U(UX))Y \ar[d]^{\Phi\otimes Y} \\
(XU)(UY) & (U(UY))X \ar[d]_\Psi \ar[rd]_{\Phi\otimes X}
& (UU)(YX) \ar[rr]^-{(UU)\otimes\Psi} \ar[d]_\Phi &
& (UU)(XY) \ar[r]^-\Phi &
((UU)X)Y \ar[dl]^{\Psi\otimes Y} \\
& X(U(UY)) \ar[lu]^\Phi
\ar[rd]_-{X\otimes\Phi} & ((UU)Y)X \ar[d]_\Psi & &
(X(UU))Y \\
& & X((UU)Y) \ar[rru]_-\Phi.
}$$
However, the leftmost and the lower triangular subdiagrams
commute by compatibility, and the upper rightmost subdiagram
commutes by coherence. The lower leftmost subdiagram commutes
by naturality of $\Psi$, and the central square sudiagram
commutes by naturality of $\Phi$. The remaining central
subdiagram commutes by coherence, and the top rectangular
subdiagram is of the form $U\otimes D$, where $D$ is a diagram
that commutes by virtue of lemma \ref{lem_second-compat}.

The uniqueness of $u_X$ is clear by inspecting the second
diagram, with $Y=U$.
\end{proof}

\begin{remark}
(i)\ \ 
Keep the notation of proposition \ref{prop_unit-isom}.
As a consequence of the naturality of $u_X$, we get a
commutative diagram
$$
\xymatrix{ X \ar[rr]^-{u_X} \ar[d]_{u_X} & &
U\otimes X \ar[d]^{u\otimes u_X} \\
U\otimes X \ar[rr]^-{u_{U\otimes X}} & & U\otimes(U\otimes X)
}$$
for every object $X$ of $\cC$. In other words :
$$
u_{U\otimes X}=U\otimes u_X
\qquad
\text{for every $X\in\Ob(\cC)$}.
$$

(ii)\ \
Let $X=Y=U$ in the second diagram of proposition
\ref{prop_unit-isom}; we obtain a commutative diagram
$$
\xymatrix{ U\otimes U \ar[rr]^-{u\otimes U} \ar[d]_{U\otimes u} & &
(U\otimes U)\otimes U \ar[d]^{\Psi_{U,U}\otimes U} \\
U\otimes(U\otimes U) \ar[rr]^-{\Phi_{U,U,U}} & &
(U\otimes U)\otimes U
}$$
Since $u=u_U$, we may combine with (i), to deduce that
$\Psi_{U,U}\otimes U=\one_{(U\otimes U)\otimes U}=
\one_{U\otimes U}\otimes U$. By naturality of $\Psi$,
it follows that $U\otimes\Psi_{U,U}=U\otimes\one_{U\otimes U}$,
and since \eqref{eq_auto-morphism} is an equivalence, we
conclude that
$$
\Psi_{U,U}=\one_{U\otimes U}.
$$
\end{remark}

\begin{example}\label{ex_stupid-tensor}
Let $\cC$ be any category with small $\Hom$-sets, in which
finite products are representable.
For every pair $(X,Y)$ of objects of $\cC$, pick an object
$X\otimes Y$ representing their product, and fix also two
projections $p_{X,Y}:X\otimes Y\to X$, $q_{X,Y}:X\otimes Y\to Y$
inducing an isomorphism of functors
$$
h_{X\otimes Y}\isom h_X\times h_Y
\quad :\quad
\phi\mapsto(p_{X,Y}\circ\phi,q_{X,Y}\circ\phi)
\qquad
\text{for every $Z\in\Ob(\cC)$ and $\phi\in h_{X\times Y}(Z)$}
$$
(notation of \eqref{subsec_yoneda}).
If $(X',Y')$ is another such pair, and $(g,h):(X,Y)\to(X',Y')$
any morphism in $\cC\times\cC$, then there exists a unique
morphism $f:X\otimes Y\to X'\otimes Y'$ such that
$$
(p_{X',Y'}\circ f,q_{X',Y'}\circ f)=(g\circ p_{X,Y},h\circ p_{X,Y})
$$
and we set $g\otimes h:=f$. These rules define a functor
$\otimes:\cC\times\cC\to\cC$. For every three objects $X,Y,Z$
there is a natural isomorphism
$X\otimes(Y\otimes Z)\isom(X\otimes Y)\otimes Z$ that yields
an associativity constraint for $\otimes$; namely, we let
$$
\Phi_{X,Y,Z}:=(p_{X,Y\otimes Z}\otimes(p_{Y,Z}\otimes q_{X,Y\otimes Z}))
\otimes(q_{Y,Z}\circ q_{X,Y\otimes Z}).
$$
Likewise, we get a commutativity constraint by setting
$$
\Psi_{X,Y}:=q_{X,Y}\otimes p_{X,Y}
\qquad
\text{for every $X,Y\in\Ob(\cC)$}.
$$
The verifications of the axioms of definition \ref{def_tensor-cat}
are lengthy but straightforward, and shall be left to the reader.
If $U$ is any final object of $\cC$ (example
\ref{ex_equalizers}(iv)), then there exists a unique
morphism $u:U\to U\otimes U$ which is easily seen to be an
isomorphism, and the pair $(U,u)$ yields a unit for $\otimes$.
In this way, any category with finite products and small
$\Hom$-sets is naturally endowed with a structure of tensor
category. Notice that, for this tensor structure on $\cC$
(and indeed, for most of the tensor categories that are
found in applications), the existence of functorial
isomorphisms $u_X$ fulfilling the conditions of proposition
\ref{prop_unit-isom}, is self-evident.
\end{example}

\begin{definition}\label{def_internal-hom}
Let $(\cC,\otimes,\Phi,\Psi)$ be a tensor category,
$X\in\Ob(\cC)$ any object, and suppose that the functor
$$
-\otimes X:\cC\to\cC
\qquad
Y\mapsto Y\otimes X
$$
admits a right adjoint :
$$
\cHom(X,-):\cC\to\cC
\qquad
Y\mapsto\cHom(X,Y).
$$
Then, we call $\cHom(X,-)$ the {\em internal $\Hom$ functor}
for the object $X$.
\end{definition}

\begin{remark}\label{rem_Hom-how-to}
(i)\ \
As usual, the internal $\Hom$ functor is determined up to
unique isomorphism, if it exists. The counit of adjunction
is a morphism of $\cC$
$$
\ev_{X,Y}:\cHom(X,Y)\otimes X\to Y
$$
called the {\em evaluation morphism}.

(ii)\ \
Suppose that every object of $\cC$ admits an internal $\Hom$
functor; then we say briefly that {\em $\cC$ admits an internal
$\Hom$ functor}, and clearly we get a functor
$$
\cC^o\times\cC\to\cC
\quad :\quad
(X,Y)\mapsto\cHom(X,Y)
\qquad
\text{for every $X,Y\in\Ob(\cC)$}.
$$
Moreover, for every $X,Y,Z\in\Ob(\cC)$ the composition
$$
\xymatrix{
(\cHom(X,Y)\otimes\cHom(Y,Z))\otimes X \ar[rr]^-\sim & &
(\cHom(X,Y)\otimes X)\otimes\cHom(Y,Z) 
\ar[d]^{\ev_{X,Y}\otimes\cHom(Y,Z)} \\
Z\xleftarrow{\ \ev_{Y,Z}\ }\cHom(Y,Z)\otimes Y & &
Y\otimes\cHom(Y,Z) \ar[ll]_-\sim
}$$
corresponds, by adjunction, to a unique {\em composition morphism}
$$
\cHom(X,Y)\otimes\cHom(Y,Z)\to\cHom(X,Z).
$$

(iii)\ \
In the situation of (ii), notice that the functor $\cC\to\cC$
given by the rule : $Z\mapsto\cHom(X,\cHom(Y,Z))$, for every
$Z\in\Ob(\cC)$, is right adjoint to the functor given by the
rule : $Z\mapsto(Z\otimes X)\otimes Y\isom Z\otimes(X\otimes Y)$.
There follows a natural isomorphism
$$
\cHom(X,\cHom(Y,Z))\isom\cHom(X\otimes Y,Z)
\qquad
\text{for every $X,Y,Z\in\Ob(\cC)$}.
$$

(iv)\ \
Moreover, for any unit $(U,u)$ of $\underline\cC$, we get
natural bijections :
$$
\Hom_\cC(U,\cHom(X,Y))\isom\Hom_\cC(U\otimes X,Y)\isom\Hom_\cC(X,Y)
\qquad
\text{for every $X,Y\in\Ob(\cC)$}.
$$
Also, for every object $Y$ of $\cC$, denote by
$u_Y:Y\isom U\otimes Y$ the isomorphism given by proposition
\ref{prop_unit-isom}; for every $X\in\Ob(X)$, it induces
natural bijections
$$
\Hom_\cC(Y,\cHom(U,X))\isom\Hom_\cC(Y\otimes U,X)
\xrightarrow{\ \Hom_\cC(u_Y\circ\Psi_{U,Y},X)\ }\Hom_\cC(Y,X)
$$
which correspond, via the Yoneda embedding, to a natural
isomorphism
$$
\cHom(U,X)\isom X
\qquad
\text{for every $X\in\Ob(\cC)$}.
$$

(v)\ \
Let $X,Y,Z$ be any three objects of $\cC$; the natural
transformation
$$
\Hom_\cC(W\otimes X,Y)\to\Hom_\cC(W\otimes(X\otimes Z),Y\otimes Z)
\qquad
\phi\mapsto(\phi\otimes Z)\circ\Phi_{W,X,Z}
$$
corresponds, via the Yoneda embedding, to a unique morphism
$$
t_{X,Y,Z}:\cHom(X,Y)\to\cHom(X\otimes Z,Y\otimes Z).
$$
The reader can check that $t_{X,Y,Z}$ also corresponds, by
adjunction, to the morphism
$$
(\ev_{X,Y}\otimes Z)\circ\Phi_{\cHom(X,Y),X,Z}:
\cHom(X,Y)\otimes(X\otimes Z)\to Y\otimes Z.
$$

(vi)\ \
In the situation of remark \ref{rem_Mac-Lane}(ii), suppose
that $\cC$ admits an internal $\Hom$ functor; then it is
easily seen that the resulting tensor category $\bFun(\cD,\cC)$
inherits as well an internal $\Hom$ functor, in the obvious
way.
\end{remark}

The formalism of tensor categories provides the language
to deal uniformly with the notions of algebras and their
modules that occur in various concrete settings.

\begin{definition}\label{def_A-B-modules}
Let $(\cC,\otimes,\Phi,\Psi)$ be a tensor category,
$A$ and $B$ any two objects of $\cC$.
\begin{enumerate}
\item
A {\em left $A$-module\/} (resp. a {\em right $B$-module}) is a
datum $(X,\mu_X)$, consisting of an object $X$ of $\cC$, and a
morphism in $\cC$ :
$$
\mu_X:A\otimes X\to X \qquad \text{(resp. $\mu_X:X\otimes B\to X$)}
$$
called the {\em scalar multiplication\/} of $X$.
\item
A {\em morphism of left $A$-modules\/} $(X,\mu_X)\to(X',\mu_{X'})$
is a morphism $f:X\to X'$ in $\cC$ which makes commute the diagram :
$$
\xymatrix{ A\otimes X \ar[r]^-{\mu_X} \ar[d]_{\one_A\otimes f} &
X \ar[d]^f \\
A \otimes X' \ar[r]^-{\mu_{X'}} & X'.}
$$
One defines likewise morphisms of right $B$-modules.
\item
An {\em $(A,B)$-bimodule\/} is a datum $(X,\mu^l_X,\mu^r_X)$ such that
$(X,\mu^l_X)$ is a left $A$-module, $(X,\mu^r_X)$ is a right
$B$-module, and the scalar multiplications commute, {\em i.e.}
the diagram
$$
\xymatrix{A\otimes(X\otimes B) \ar[rrrr]^-{\one_A\otimes\mu^r_X}
\ar[d]_{\Phi_{A,X,B}} & & & & A\otimes X \ar[d]^{\mu^l_X} \\
(A\otimes X)\otimes B \ar[rr]^-{\mu^l_X\otimes\one_B} & &
X\otimes B \ar[rr]^-{\mu^r_X} & & X
}$$
commutes. Of course, a morphism of $(A,B)$-bimodules must be
compatible with both left and right multiplication.
\end{enumerate}
We denote by $A\Mod_l$ (resp. $B\Mod_r$, resp. $(A,B)\Mod$)
the category of left $A$-modules (resp. right $B$-modules, resp.
$(A,B)$-bimodules). For any two left $A$-modules (resp. right
$B$-modules, resp. $(A,B)$-bimodules) $X$ and $X'$, we shall
write
$$
\Hom_{A_l}(X,X')
\qquad
\text{(\ resp.\ $\Hom_{B_r}(X,X')$\ )}
\qquad
\text{(\ resp.\ $\Hom_{(A,B)}(X,X')$\ )}
$$
for the set of morphisms of left $A$-modules (resp. of right
$B$-modules, resp. of $(A,B)$-bimodules) $X\to X'$.
\end{definition}

\sset\subsubsection{}\label{subsec_Hom-as-equal}
In the situation of definition \ref{def_A-B-modules}, notice that
$$
\Hom_{B_r}(X,X')=\Equal(
\xymatrix{\Hom_\cC(X,X') \ar@<.5ex>[r]^-\alpha \ar@<-.5ex>[r]_-\beta
& \Hom_\cC(X\otimes B,X')})
$$
where $\alpha$ (resp. $\beta$) is given by the rule :
$$
f\mapsto f\circ\mu_X
\qquad
\text{(\ resp. $f\mapsto\mu_{X'}\circ(f\otimes B)$\ )}
\qquad
\text{for every $f\in\Hom_\cC(X,X')$}
$$
and similarly for left $A$-modules.
Now, suppose that all equalizers in $\cC$ are representable,
and that $\cC$ admits an internal $\Hom$ functor; then we may
define
$$
\cHom_{B_r}(X,X'):=\Equal(
\xymatrix{\cHom(X,X') \ar@<.5ex>[r]^-\alpha \ar@<-.5ex>[r]_-\beta
& \cHom(X\otimes B,X')})
$$
where $\alpha:=\cHom(\mu_X,X')$ and
$\beta:=\cHom(X\otimes B,\mu_{X'})\circ t_{X,X',B}$
(notation of remark \ref{rem_Hom-how-to}(v)). Then, it is easily
seen that the bijections of remark \ref{rem_Hom-how-to}(iv)
induce natural identifications
$$
\Hom_\cC(U,\cHom_{B_r}(X,X'))\isom\Hom_{B_r}(X,X')
\qquad
\text{for every $X,X'\in\Ob(B_r\Mod)$}.
$$
Likewise we may represent in $\cC$ the set of morphisms
between two left $A$-modules, and two $(A,B)$-modules
(details left to the reader).

\sset\subsubsection{}\label{subsec_ABC-mod-Hom}
Let $\underline\cC$ be a tensor category as in
\eqref{subsec_Hom-as-equal} and $A,B,C$ any three objects
of $\cC$; suppose that $(X,\mu_X^l,\mu_X^r)$ is an
$(A,B)$-bimodule, and $(X',\mu^l_{X'},\mu^r_{X'})$ a
$(C,B)$-bimodule. Then we claim that
$\cH:=\cHom_{B_r}((X,\mu^r_X),(X',\mu^r_{X'}))$ is naturally
a $(C,A)$-bimodule. For this, we have to exhibit natural
morphisms
$$
C\otimes\cH\xrightarrow{\ \mu_l\ }\cH\xleftarrow{\ \mu_r\ }\cH\otimes A
$$
fulfilling the condition of definition \ref{def_A-B-modules}(iii).
However, by adjunction, the datum of $\mu_l$ is the same as that
of a morphism $C\to\cHom(\cH,\cH)$, and since the functor
$\cHom(X,-)$ is left exact, the latter is the same as a
morphism
$$
C\to\Equal(\xymatrix{\cHom(\cH,\cHom(X,X'))
\ar@<.5ex>[rr]^-{\cHom(\cH,\alpha)}
\ar@<-.5ex>[rr]_-{\cHom(\cH,\beta)} & &
\cHom(\cH,\cHom(X\otimes B,X'))
}$$
which in turn -- by remark \ref{rem_Hom-how-to}(iii) --
corresponds to a morphism
$$
C\to\Equal(\xymatrix{\cHom(\cH\otimes X,X')
\ar@<.5ex>[rrr]^-{\cHom(\cH\otimes\alpha,X')}
\ar@<-.5ex>[rrr]_-{\cHom(\cH\otimes\beta,X')} & & &
\cHom(\cH\otimes(X\otimes B),X')
}$$
and again, the latter is the same as an element of
$$
\Equal(\xymatrix@C-4pt{\Hom_\cC(C\otimes(\cH\otimes X),X')
\ar@<.5ex>[rrrr]^-{\Hom_\cC(C\otimes(\cH\otimes\alpha),X')}
\ar@<-.5ex>[rrrr]_-{\Hom_\cC(C\otimes(\cH\otimes\beta),X')} & & & &
\Hom_\cC(C\otimes(\cH\otimes(X\otimes B)),X').
}$$
By unwinding the definition, it is easily seen that
the composition
$$
\bar\mu_l:
C\otimes(\cH\otimes X)\xrightarrow{\ \ev_{X,X'}\ }
C\otimes X'\xrightarrow{\ \mu^l_{X'}\ }X'
$$
lies in the above equalizer, and it provides a left
$C$-module structure for $\cH$. Likewise, $\mu_r$ shall
be the morphism corresponding to the composition
$$
\bar\mu_r:
(\cH\otimes X)\otimes A\isom\cH\otimes(A\otimes X)
\xrightarrow{\ \cH\otimes\mu_X^l\ }\cH\otimes X
\xrightarrow{\ \ev_{X,X'}\ }X'.
$$
Then, the condition that $(\cH,\mu_l,\mu_r)$ is a
bimodule, comes down to the commutativity of the
diagram
$$
\xymatrix{C\otimes((\cH\otimes X)\otimes A)
\ar[rrrrr]^-{\one_C\otimes\bar\mu_r}
\ar[d]_{C\otimes\Phi_{\cH,X,A}} & & & & & C\otimes X' \ar[d]^{\mu^l_{X'}} \\
C\otimes(\cH\otimes(X\otimes A))
\ar[rrrr]^-{C\otimes(\cH\otimes(\mu^l_X\circ\Psi_{X,A}))} & & & &
C\otimes(\cH\otimes X) \ar[r]^-{\bar\mu_l} & X'
}$$
which is immediate (details left to the reader).
We have thus obtained a bifunctor :
\set\begin{equation}\label{eq_hom-fctrs}
\cHom_{B_r}(-,-):(A,B)\Mod^o\times(C,B)\Mod\to(C,A)\Mod.
\end{equation}
Likewise, we may define a bifunctor :
$$
\cHom_{A_l}(-,-):(A,B)\Mod^o\times(A,C)\Mod\to(B,C)\Mod.
$$

\sset\subsubsection{}\label{subsec_adjunct-scoppia}
Keep the situation of \eqref{subsec_ABC-mod-Hom}, and
suppose moreover that all coequalizers in $\cC$ are
representable. Fix an $(A,B)$-bimodule
$(X,\mu_X^l,\mu_X^r)$; the functor
$$
(C,B)\Mod\to(C,A)\Mod
\quad : \quad
X'\mapsto\cHom_{B_r}(X,X')
$$
admits a left adjoint, the {\em tensor product}
$$
(C,A)\Mod\to(C,B)\Mod
\quad :\quad
(X',\mu^l_{X'},\mu^l_{X'})\mapsto
(X',\mu^l_{X'},\mu^r_{X'})\otimes_A(X,\mu_X^l,\mu_X^r)
$$
given by the coequalizer (in $\cC$) of the morphisms :
$$
\xymatrix{ X'\otimes(A\otimes X)
\ar@<.5ex>[rrr]^-{\one_{X'}\otimes\mu^l_X}
\ar@<-.5ex>[rrr]_-{(\mu^r_{X'}\otimes\one_X)\circ\Phi_{X',A,X}}
& & & X'\otimes X}
$$
with scalar multiplications induced by $\mu^r_X$ and $\mu^l_{X'}$.
Likewise, we have a functor :
$$
(A,C)\Mod\to(B,C)\Mod \quad : \quad
(X',\mu^l_{X'},\mu^l_{X'})\mapsto
(X,\mu_X^l,\mu_X^r)\otimes_A(X',\mu^l_{X'},\mu^l_{X'})
$$
which admits a similar description, and is left adjoint to
the functor $X'\mapsto\cHom_{A_l}(X,X')$ (verifications left
to the reader).

\sset\subsubsection{}\label{subsec_adjunct-scoppia-right}
Let $(U,u)$ be unit object for $\underline\cC$. Notice that,
for any object $A$ of $\cC$, the rule
$(Y,\mu_Y)\mapsto(Y,u^{-1}_Y,\mu_Y)$ (where $u_Y:Y\to U\otimes Y$
is the natural isomorphism supplied by proposition
\ref{prop_unit-isom}), induces a faithful functor
$A\Mod_r\to(U,A)\Mod$. Letting $C:=U$ in
\eqref{subsec_adjunct-scoppia}, we see that any $(A,B)$-bimodule
$X$ also determines a functor :
$$
A\Mod_r\to B\Mod_r
\quad : \quad
Y\mapsto Y\otimes_AX
$$
and likewise for left $A$-modules.

\begin{example}\label{ex_cp}
If $\cC=\Set$ is the category of sets (regarded as a tensor category
as in example \ref{ex_stupid-prior}), then a left $A$-module is just
a set $X'$ with a {\em left action} of $A$, {\em i.e.} a map of sets
$$
A\times X'\to X'
\quad : \quad
(a,x)\mapsto a\cdot x.
$$
An $(A,A)$-bimodule $X$ is a set with both left and right actions
of $A$, such that $(a\cdot x)\cdot a'=a\cdot(x\cdot a')$ for every
$a,a'\in A$ and every $x\in X$. With this notation, the tensor
product $X'\otimes_AX$ is the quotient $(X'\times X)\!/\!\!\sim$,
where $\sim$ is the smallest equivalence relation such that
$(x'a,x)\sim(x',ax)$ for every $x\in X$, $x'\in X'$ and $a\in A$.
\end{example}

\begin{definition}\label{def_T-not-only-top}
Let $\underline\cC:=(\cC,\otimes,\Phi,\Psi)$ be a tensor
category, and  $(U,u)$ a unit for $\underline\cC$. 
\begin{enumerate}
\item
A {\em $\underline\cC$-semigroup} is a datum $(M,\mu_M)$
consisting of an object $M$ of $\cC$ and a morphism
$\mu_M:M\otimes M\to M$, the {\em multiplication law\/} of
$\underline M$, such that $(M,\mu_M,\mu_M)$ is a $(M,M)$-bimodule.
A morphism of $\underline\cC$-semigroups is a morphism
$\phi:M\to M'$ in $\cC$, such that
$$
\mu_{M'}\circ(\phi\otimes\phi)=\phi\circ\mu_M.
$$
\item
A {\em $\underline\cC$-monoid\/} is a datum
$\underline M:=(M,\mu_M,1_M)$, where $(M,\mu_M)$ is a
semigroup, and $1_M:U\to M$ is a morphism in $\cC$, called
the {\em unit\/} of $\underline M$, such that
$$
\mu_M\circ(1_M\otimes\one_M)\circ u_M=\one_M=
\mu_M\circ(\one_M\otimes 1_M)\circ u_M
$$
where $u_M:M\isom U\otimes M$ is the natural isomorphism
provided by proposition \ref{prop_unit-isom}.
We say that $\underline M$ is {\em commutative}, if
$$
\mu_M=\mu_M\circ\Psi_{M,M}.
$$
A morphism of monoids $\underline M\to\underline M'$ is a
morphism of semigroups $\phi:M\to M'$ such that
$\phi\circ 1_M=1_{M'}$.
\end{enumerate}
\end{definition}

\begin{example}\label{ex_stupid-prior}
(i)\ \ 
Let $\cC$ be any category with small $\Hom$-sets,
in which finite products are representable, endow $\cC$
with the tensor category structure described in example
\ref{ex_stupid-tensor}, and pick a final object $1_\cC$
of $\cC$. Then, a $\underline\cC$-monoid is a datum
$\underline M:=(M,\mu_M,1_M)$, where $\mu_M:M\times M\to M$
and $1_M:1_\cC\to M$ are morphisms of $\cC$, and the
axioms for $\mu_M$ and $1_M$ can be rephrased as requiring
that, for every object $X$ of $\cC$, the set
$M(X):=\Hom_\cC(X,M)$, endowed with the composition law :
$$
M(X)\times M(X)\isom\Hom_\cC(X,M\times M)
\xrightarrow{\ \Hom_\cC(X,\mu_M)\ }
M(X)
\qquad
(m,m')\mapsto m\cdot m'
$$
is a (usual) monoid, with unit $\Img\,1_M(X)\in M(X)$.
Of course, $\underline M$ is commutative, if and only if
$m\cdot m'=m'\cdot m$ for all objects $X$ of $\cC$ and every
$m,m'\in M(X)$.

(ii)\ \
In the situation of (i), a $\underline\cC$-monoid shall also be
called simply a $\cC$-monoid. The category of $\cC$-monoids
admits an initial object which is also a final object, namely
$\underline 1{}_\cC:=(1_\cC,\mu_1,\one_1)$, where $\mu_1$ is
the (unique) morphism $1_\cC\times 1_\cC\to 1_\cC$.
(Most of the above can be repeated with the theory of semigroups
replaced by any "algebraic theory" in the sense of
\cite[Def.3.3.1]{BorII} : {\em e.g.} in this way one can define
$\cC$-groups, $\cC$-rings, and so on.)
\end{example}

\sset\subsubsection{}\label{subsec_sheaves-of-mons}
Let $\underline\cC$ and $U$ be as in definition
\ref{def_T-not-only-top}, and $\underline M:=(M,\mu_M,1_M)$
a $\underline\cC$-monoid; of course, we are especially
interested in the $M$-modules which are compatible with
the unit and multiplication law of $M$. Hence we define
a {\em left $\underline M$-module} as a left $M$-module
$(S,\mu_S)$ such that the following diagrams commute :
$$
\xymatrix{
U\otimes S \ar@{=}[r] \ar[d]_{1_M\otimes\one_S} & U\otimes S 
& M\otimes(M\otimes S) \ar[rrr]^-{\one_M\otimes\mu_S}
\ar[d]_{\Phi_{M,M,S}}
& & & M\otimes S \ar[d]^{\mu_S} \\
M\otimes S \ar[r]^-{\mu_S} & S \ar[u]_{u_S} &
(M\otimes M)\otimes S \ar[rr]^-{\mu_M\otimes\one_S} & &
M\otimes S \ar[r]^-{\mu_S} & S
}$$
where $u_S$ is the isomorphism given by proposition
\ref{prop_unit-isom}.
Likewise we define right $\underline M$-modules, and
$(\underline M,\underline N)$-bimodules, if $\underline N$
is a second $\underline\cC$-monoid; especially,
$(\underline M,\underline M)$-bimodules shall also be called
simply $\underline M$-bimodules.

A morphism of left $\underline M$-modules 
$(S,\mu_S)\to(S',\mu_{S'})$ is just a morphism of left
$M$-modules, and likewise for right modules and bimodules.
For instance, $\underline M$ is a $\underline M$-bimodule in a
natural way, and an {\em ideal\/} of $M$ is a
sub-$\underline M$-bimodule $I$ of $\underline M$.
We denote by $\underline M\Mod_l$ (resp. $\underline M\Mod_r$,
resp. $\underline M\Mod$) the category of left (resp. right,
resp. bi-) $\underline M$-modules; more generally, if $\underline N$
is a second $\underline\cC$-monoid, we have the category
$(\underline M,\underline N)\Mod$ of the corresponding bimodules.

\begin{example}\label{ex_prior}
Take $\cC:=\Set$, regarded as a tensor category, as in
example \ref{ex_stupid-tensor}. Then a $\cC$-monoid is
just a usual monoid $M$, and a left $M$-module is a
datum $(S,\mu_S)$ consisting of a set $S$ and a
{\em scalar multiplication\/}
$M\times S\to S$ : $(m,s)\mapsto m\cdot s$ such that
$$
1\cdot s=s\quad
\text{and}\quad
x\cdot(y\cdot s)=(x\cdot y)\cdot s
\qquad
\text{for every $x,y\in M$ and every $s\in S$}.
$$
A morphism $\phi:(S,\mu_S)\to(T,\mu_T)$ of $M$-modules
is then a map of sets $S\to T$ such that
$$
x\cdot\phi(s)=\phi(x\cdot s)
\qquad
\text{for every $x\in M$ and every $s\in S$}
$$
Likewise, an ideal of $M$ is a subset $I\subset M$ such that
$a\cdot x, x\cdot a\in I$ whenever $a\in I$ and $x\in M$.
\end{example}

\begin{remark}\label{rem_commutes-forgets}
Let $\cC$ be a complete and cocomplete category, whose colimits
are universal (see example \ref{ex_universal-col}), and
$\underline M$ a $\cC$-monoid (see example
\ref{ex_stupid-prior}(ii)).

(i)\ \
The categories $\underline M\Mod_l$, $\underline M\Mod_r$ and
$(\underline M,\underline N)\Mod$ are complete and cocomplete,
and the forgetful functor $\underline M\Mod_l\to\cC$ (resp. the
same for right modules and bimodules) commutes with all limits
and colimits.

(ii)\ \
Notice also that the forgetful functor
$\underline M\Mod_l\to\cC$ is conservative. Together with (i)
and propositions \ref{prop_was-also-cofinal} and
\ref{prop_conser-nd-left-exact}(i), this implies
that a morphism of left $\underline M$-modules is a monomorphism
(resp. an epimorphism) if and only if the same holds for the
underlying morphism in $\cC$ (and likewise for right modules
and bimodules).

(iii)\ \
For each of these categories, the initial object is just
the initial object $\emptyset_\cC$ of $\cC$, endowed with
the trivial scalar multiplication.
Likewise, the final object is the final object $1_\cC$ of $\cC$,
with scalar multiplication given by the unique morphism
$M\times 1_\cC\to 1_\cC$. Moreover, the forgetful functor
$\underline M\Mod_l\to\cC$ admits a left adjoint, that
assigns to any $\Sigma\in\Ob(\cC)$ the {\em free
$\underline M$-module\/} $\underline M^{(\Sigma)}$
generated by $\Sigma$; as an object of $\cC$, the latter
is just $M\times\Sigma$, and the scalar multiplication
is derived from the composition law of $\underline M$,
in the obvious way.

For instance, for any $n\in\N$, and any left (or right or bi-)
$\underline M$-module $S$, we denote as usual by $S^{\oplus n}$
the coproduct of $n$ copies of $S$.
\end{remark}

\begin{remark}
Let $\underline\cC$ be a tensor category, $\underline M$,
$\underline N$, $\underline P$, $\underline Q$ four
$\underline\cC$-monoids.

(i)\ \
Let $S$ be a $(\underline M,\underline N)$-bimodule, $S'$ a
$(\underline P,\underline N)$-bimodule and $S''$ a
$(\underline P,\underline M)$-bimodule. Then it is easily
seen that the $(P,M)$-bimodule (resp. the $(P,N)$-bimodule)
$\cHom_{N_r}(S,S')$ (resp. $S''\otimes_MS$) is actually a
$(\underline P,\underline M)$-bimodule (resp. a
$(\underline P,\underline N)$-bimodule) and the adjunction of
\eqref{subsec_adjunct-scoppia} restricts to an adjunction between
the corresponding categories of bimodules : the details shall be
left to the reader.

(ii)\ \
We have as well the analogue of the usual associativity constraints.
Namely, for every $(\underline M,\underline N)$-bimodule $S$, every
$(\underline N,\underline P)$-bimodule $S'$ and every
$(\underline P,\underline Q)$-bimodule $S''$, there is a natural
isomorphism
$$
(S\otimes_NS')\otimes_PS''\isom S\otimes_N(S'\otimes_PS'')
\qquad
\text{in $(\underline M,\underline Q)\Mod$}
$$
and natural isomorphisms $M\otimes_MS\isom S\isom S\otimes_NN$
in $(\underline M,\underline N)\Mod$.

(iii)\ \
Also, if $\underline M$ is commutative, every left (or right)
$\underline M$-module is naturally a
$(\underline M,\underline M)$-bimodule, and we have a commutative
constraint
$$
S\otimes_MS'\isom S'\otimes_MS
\qquad
\text{for all left (or right) $\underline M$-modules}.
$$
And taking into account (ii), it is easily seen that
$(\underline M\Mod_l,\otimes_M)$ is a tensor category.
\end{remark}

\sset\subsubsection{}\label{subsec_tens-restr}
Let $\phi:\underline M_1\to\underline M_2$ be a morphism of
$\underline\cC$-monoids; we have the
$(\underline M_1,\underline M_2)$-bimodule :
$$
M_{1,2}:=(M_2,\mu_{M_2}\circ(\phi\otimes\one_{M_2}),\mu_{M_2}).
$$
Letting $X:=M_{1,2}$ in \eqref{subsec_adjunct-scoppia-right},
we obtain a {\em base change functor\/} for right modules :
$$
\underline M_1\Mod_r\to\underline M_2\Mod_r
\quad : \quad
X\mapsto X\otimes_{M_1}M_2:=X\otimes_{M_1}M_{1,2}.
$$
The base change is left adjoint to the {\em restrictions of
scalars\/} associated with $\phi$, {\em i.e.} the functor :
$$
\underline M_2\Mod_r\to\underline M_1\Mod_l
\quad : \quad
(X,\mu_X)\mapsto(X,\mu_X)_{(\phi)}:=(X,\mu_X\circ(\one_X\otimes\phi))
$$
(verifications left to the reader). The same can be repeated, as
usual, for left modules; for bimodules, one must take the tensor
product on both sides : $X\mapsto M_2\otimes_{M_1}X\otimes_{M_1}M_2$.

\begin{example}\label{ex_rank-mod}
Take $\cC=\Set$, and let $M$ be any monoid, $\Sigma$ any set,
and $M^{(\Sigma)}$ the free $M$-module generated by $\Sigma$.
From the isomorphism
$$
M^{(\Sigma)}\otimes_M\{1\}\isom\{1\}^{(\Sigma)}=\Sigma
$$
we see that the cardinality of $\Sigma$ is an invariant, called the
{\em rank} of the free $M$-module $M^{(\Sigma)}$, and which we
denote $\rk_MM^{(\Sigma)}$.
\end{example}

\begin{definition}
(i)\ \ 
A {\em pre-additive category} is the datum of a category $\cA$
and of an abelian group structure on $\Hom_\cA(A,B)$ for every
$A,B\in\Ob(\cA)$ (especially, the $\Hom$-sets of $\cA$ are not
empty), such that the following holds. For every
$A,B,C\in\Ob(\cA)$, the composition law
$$
\Hom_\cA(A,B)\times\Hom_\cA(B,C)\to\Hom_\cA(A,C)
$$
is a bilinear pairing.

(ii)\ \
A functor $F:\cA\to\cB$ between pre-additive categories is
{\em additive} if it induces group homomorphisms
$$
\Hom_\cA(X,Y)\to\Hom_\cB(FX,FY)
\quad : \quad
\phi\mapsto F\phi
\qquad
\text{for every $X,Y\in\Ob(\cA)$}.
$$
\end{definition}

\begin{remark}\label{rem_additive-cat}
Let $\cA$ be any pre-additive category, and choose a
universe $\sU$ such that $\cA$ is $\sU$-small, so that
all the finite limits and finite colimits of $\cA$ are
well defined as $\sU$-presheaves.

(i)\ \
If $A\in\Ob(\cA)$ is any object, denote by $\zero_A$ the
neutral element of the abelian group $\End_\cA(A)$. Suppose
that the equalizer of the pair of morphisms
$\one_A,\zero_A:A\to A$ is representable by an object $0$
of $\cA$ (see example \ref{ex_equalizers}(ii)). Then, the
datum of a morphism $B\to 0$ is the same as that of a
morphism $\phi:B\to A$ that factors through $\zero_A$.
By the bilinearity of the $\Hom$-pairing, the latter condition
holds if and only if $\phi$ is the neutral element of
$\Hom_\cA(B,A)$. We conclude that $0$ is a final object
in $\cA$. Dually, if the coequalizer of the pair
$(\one_A,\zero_A)$ is representable by some object $0'$
of $\cA$, then $0'$ is initial in $\cA$.
Moreover, if $\cA$ admits a final object $0$, then it
is easily seen that the unique morphism $A\to 0$ is also
the coequalizer of the pair $(\one_A,\zero_A)$, so $0$ is
also an initial object. Conversely, if $\cA$ admits an
initial object, then this object is also final in $\cA$,
and for any two objects $A,B$ of $\cA$, the neutral
element $\zero_{A,B}$ of $\Hom_\cA(A,B)$ is the unique
morphism that factors through $0$. We say that $0$ is
a {\em zero object} for $\cA$.

(ii)\ \
Suppose that the product $A_1\times A_2$ is representable in
$\cA$ for given $A_1,A_2\in\Ob(\cA)$. Denote by
$p_i:A_1\times A_2\to A_i$ ($i=1,2$) the projections; then,
there are unique morphisms $e_i:A_i\to A_1\times A_2$ ($i=1,2$)
such that
\set\begin{equation}\label{eq_proj-inj-maps}
p_i\circ e_i=\one_{A_i}
\qquad
\text{for $i=1,2$}
\qquad\text{and}\qquad
p_i\circ e_j=\zero_{A_j,A_i}
\qquad
\text{for $i\neq j$}.
\end{equation}
Notice that
\set\begin{equation}\label{eq_prod-is-coprod}
e_1\circ p_1+e_2\circ p_2=\one_{A_1\times A_1}.
\end{equation}
Indeed, we have
$$
p_i\circ(e_1\circ p_1+e_2\circ p_2)=(p_i\circ e_1\circ p_1)+
(p_i\circ e_2\circ p_2)=p_i
\qquad
\text{$i=1,2$}
$$
by bilinearity of the $\Hom$ pairing, and $\one_{A_1\times A_2}$
is the unique endomorphism $\phi$ of $A_1\times A_2$ such that
$p_i\circ\phi=p_i$ for $i=1,2$. It follows that $A_1\times A_2$
also represents the coproduct $A_1\amalg A_2$. Indeed, say that
$f_i:A_i\to B$, for $i=1,2$, are two morphisms to another object
$B$ of $\cA$, and set
$f:=f_1\circ p_1+f_2\circ p_2:A_1\times A_2\to B$; it is easily
seen that $f\circ e_i=f_i$ for $i=1,2$, and by virtue of
\eqref{eq_prod-is-coprod}, the morphism $f$ is the unique
one that satisfies these identities. Conversely, if the
coproduct of $A_1$ and $A_2$ is representable, a similar
argument shows that also $A_1\times A_2$ is representable.
We say that $A_1\times A_2$ is a {\em biproduct} of $A_1$
and $A_2$, and denote it by $A_1\oplus A_2$.

(iii)
Notice that the morphisms $(p_i,e_i~|~i=1,2)$ with the
identities \eqref{eq_proj-inj-maps} and \eqref{eq_prod-is-coprod}
characterize $A_1\oplus A_2$ up to unique isomorphism.
Namely, say that $B$ is another object of $\cA$, for
which exist morphisms $p'_i:B\to A_i$ and $e'_i:A_i\to B$
($i=1,2$) such that $p'_i\circ e'_i=\one_{A_i}$ for $i=1,2$,
and $p'_i\circ e'_j=\zero_{A_jA_i}$ for $i\neq j$, and
moreover $e'_1\circ p'_1+e'_2\circ p'_2=\one_B$. Then
the pair $(e'_1,e'_2)$ (resp. $(p'_1,p'_2)$) induces a
morphism $e':A_1\oplus A_2\to B$ (resp. $p':B\to A_1\oplus A_2$)
with
$$
p_j\circ p'\circ e'\circ e_i=p'_j\circ e'_i=p_j\circ e_i
\qquad
\text{for $i,j=1,2$}
$$
which -- by virtue of the universal properties of the biproduct --
implies that $p'\circ e'=\one_{A_1\oplus A_2}$. Likewise, we
may compute
$$
\begin{aligned}
e'\circ p'=\, &
(e'\circ e_1\circ p_1+e'\circ e_2\circ p_2)\circ
(e_1\circ p_1\circ p'+e_2\circ p_2\circ p') \\
=\, &
e'\circ e_1\circ p_1\circ p'+e'\circ e_2\circ p_2\circ p' \\
=\, &
e'_1\circ p'_1+e'_2\circ p'_2=\one_B
\end{aligned}
$$
whence the contention.

(iv)\ \
Suppose that $B_1$ and $B_2$ are any other two objects of
$\cA$ such that $B_1\oplus B_2$ is also representable;
given two morphisms $f_1:A_1\to B_1$ and $f_2:A_2\to B_2$,
we denote by $f_1\oplus f_2:A_1\oplus A_2\to B_1\oplus B_2$
the unique morphism such that
$$
p_{B,i}\circ(f_1\oplus f_2)\circ e_{A,i}=f_i
\quad
\text{for $i=1,2$, and}
\quad
p_{B,i}\circ(f_1\oplus f_2)\circ e_{A,j}=\zero_{A_j,B_i}
\quad
\text{for $i\neq j$}.
$$
Notice that
\set\begin{equation}\label{eq_sum-is-biprod}
f_1\oplus f_2=(f_1\oplus\zero_{A_2,B_2})+(\zero_{A_1,B_1}\oplus f_2)
\end{equation}
(where the sum is taken in the abelian group
$\Hom_\cA(A_1\oplus A_2,B_1\oplus B_2)$); indeed, by bilinearity
of the $\Hom$ pairing, it is easily seen that the right-hand
side of \eqref{eq_sum-is-biprod} also satisfies the same
identities above that define $f_1\oplus f_2$.

(v)\ \
If $f:A\to B$ is any morphism of $\cA$, then we define
the {\em kernel} (resp. {\em cokernel}) of $f$ as the
equalizer (resp. coequalizer)
$$
\Ker\,f:=\Equal(f,\zero_{A,B})
\qquad
\Coker\,f:=\mathrm{Coequal}(f,\zero_{A,B}).
$$
Suppose that all kernels and cokernels of $\cA$ are representable
in $\cA$, and denote by
$$
\iota_f:\Ker\,f\to A
\qquad
\pi_f:B\to\Coker\,f
$$
the natural morphisms in $\cA$. Notice that $\iota_f$ is a monomorphism,
and $\pi_f$ an epimorphism. Notice also that $f$ factors uniquely
in $\cA$ as a composition
\set\begin{equation}\label{eq_ker-coker}
A\xrightarrow{\ \pi_{\iota_f}\ }\Coker\,\iota_f
\xrightarrow{\ \beta_f\ }\Ker\,\pi_f
\xrightarrow{\ \iota_{\pi_f}\ } B.
\end{equation}

(vi)\ \
Suppose that $\cA$ admits a zero object $0$, and let $f:A\to B$
be any morphism; by definition $\Ker\,f$ is the presheaf such that
$$
\Ker\,f(C)=\{g:B\to C~|~g\circ f=g\circ\zero_{A,B}=\zero_{A,C}\}.
$$
If $f$ is a monomorphism, the identity
$g\circ f=\zero_{A,C}=\zero_{B,C}\circ f$ implies that
$g=\zero_{B,C}$, so $\Ker\,f$ is represented by $0$.
Dually, if $f$ is an epimorphism, then $\Coker\,f$ is
represented by $0$.
\end{remark}

\begin{remark}\label{rem_fun-exact-add}
Let $\cA$, $\cB$ be any two pre-additive categories that
admit a zero object, and $F:\cA\to\cB$ a functor.

(i)\ \ 
If $F$ is additive, remark \ref{rem_additive-cat}(iii)
immediately implies that $F$ transforms representable
biproducts into representable biproducts. The latter
assertion still holds in case $F$ is not necessarily
additive, but is either left or right exact. Indeed,
suppose that $F$ is left exact, let $A_1\oplus A_2$
be any biproduct, and let $p_i,e_i$ be the morphisms
described in remark \ref{rem_additive-cat}(ii); by left
exactness, $F(A_1\oplus A_2)$ represents the product of
$FA_1$ and $FA_2$, whence an isomorphism
$FA_1\oplus FA_2\isom F(A_1\oplus A_2)$ which identifies
$Fp_1$ and $Fp_2$ with the natural projections.
Moreover, $F$ transforms the final object of $\cA$ into
the final object of $\cB$ (see example
\ref{ex_equalizers}(iv)); then, by inspecting
the argument in remark \ref{rem_additive-cat}(ii), it is
easily seen that $F$ identifies as well $Fe_i$ with the
natural injections $FA_i\to FA_1\oplus FA_2$, for $i=1,2$,
so the assertion follows from remark \ref{rem_additive-cat}(iii).
A similar argument works in case $F$ is right exact.

(ii)\ \
Suppose moreover, that all biproducts of $\cA$ are representable.
Then we claim that the abelian group structure on $\Hom_\cA(A,B)$
is determined by the category $\cA$, {\em i.e.} if $\cB$ is
any other pre-additive category, and $F:\cA\to\cB$ is any
equivalence of categories, then $F$ induces group isomorphisms
(and not just bijections) on $\Hom$ sets. Indeed, let $A$ and
$B$ be any two objects of $\cA$, and denote by
$\Delta_A:A\to A\oplus A$ (resp. $\mu_B:B\oplus B\to B$) the
unique morphism such that $p_i\circ\Delta_A=\one_A$ (resp.
$\mu_B\circ e_i=\one_B$) for $i=1,2$. Then we have
$$
f_1+f_2=\mu_B\circ(f_1\oplus f_2)\circ\Delta_A
\qquad
\text{for every $f_1,f_2:A\to B$}
$$
where $f_1+f_2$ denotes the sum in the abelian group
$\Hom_\cA(A,B)$. Indeed, since clearly
$\zero_{A_1,B_1}\oplus\zero_{A_2,B_2}=\zero_{A_1\oplus A_2,B_1\oplus B_2}$,
identity \eqref{eq_sum-is-biprod} reduces to checking
that $f_1=\mu_B\circ(f_1\oplus\zero_{A_2,B_2})\circ\Delta_A$
(and likewise for $f_2$), which follows easily from
\eqref{eq_prod-is-coprod} : details left to the reader.

(iii)\ \
Combining (i) and (ii) we see that if all biproducts of
$\cA$ are representable, and $F$ is either left or right
exact, then $F$ is additive. More generally, we see that
for $F$ to be additive, it suffices that $F$ sends the
zero object of $\cA$ to the zero object of $\cB$, and that
$F$ commutes with the biproducts of the form $A\oplus A$,
for every $A\in\Ob(\cA)$.
\end{remark}

\begin{definition}
(i)\ \
An {\em additive category} is a pre-additive category
which admits a zero object, and whose biproducts are
representable.

(ii)\ \
An {\em abelian category} is an additive
category $\cA$ such that the following holds:
\begin{enumerate}
\alphaenu
\item
All the kernels and cokernels of $\cA$ are representable.
\item
For every morphism $f$ of $\cA$, the morphism $\beta_f$
of \eqref{eq_ker-coker} is an isomorphism.
\end{enumerate}
\end{definition}

\begin{example}\label{ex_idemp-in-additive-cat}
Let $\cA$ be an additive category, $A$ an object of $\cA$,
and $e:A\to A$ an endomorphism of $A$ such that :
\begin{enumerate}
\alphaenu
\item
$e=e\circ e$, {\em i.e.} $e$ is an {\em idempotent} element
of the ring $\End_\cA(A)$.
\item
$\Ker\,(e)$ and $\Ker\,(\one_A-e)$ are representable in $\cA$.
\end{enumerate}
Then the morphisms $\iota_e:\Ker\,(e)\to A$ and
$\iota_{\one_A-e}:\Ker\,(\one_A-e)\to A$ induce an isomorphism
$$
\omega:\Ker\,(e)\oplus\Ker\,(\one_A-e)\isom A.
$$
Indeed, it is easily seen that the morphism $\eps:A\to A\oplus A$
defined by the pair $(\one_A-e,e)$ factors as a composition
$$
A\xrightarrow{\ \bar\eps\ }\Ker\,(e)\oplus\Ker\,(\one_A-e)
\xrightarrow{\ \iota_e\oplus\iota_{\one_A-e}}A\oplus A
$$
for a unique morphism $\bar\eps$ of $\cA$. Then a direct
computation shows that $\omega\circ\bar\eps=\one_A$ and
$\bar\eps\circ\omega=\one_{\Ker\,(e)\oplus\Ker\,(\one_A-e)}$, whence
the assertion (details left to the reader).
\end{example}

\begin{remark}\label{rem_Add.Fun}
Let $\cA$ be any pre-additive category.

(i)\ \
The category $\cA^o$ is naturally pre-additive : namely,
if $A,B\in\Ob(\cA^o)$ are any two objects, one endows
$\Hom_{\cA^o}(A,B)$ with the group structure of
$\Hom_\cA(B,A)$. Likewise, if $F:\cA\to\cB$ is any
additive functor between the pre-additive categories
$\cA$, $\cB$, then also $F^o:\cA^o\to\cB^o$ is additive.
Moreover, if $\cA$ is additive the same holds for $\cA^o$.
Indeed, the zero object of $\cA$ is obviously a zero
object also in $\cA^o$, and biproducts are representable
in $\cA^o$, since one can take
$$
A^o\oplus B^o:=(A\oplus B)^o
\qquad
\text{for every $A,B\in\Ob(\cA)$}.
$$
Furthermore, if $\cA$ is abelian, the same holds for $\cA^o$.
Indeed, we can take
$$
\Ker\,f^o:=(\Coker\,f)^o
\qquad
\Coker\,f^o:=(\Ker\,f)^o
\qquad
\text{for every morphism $f$ in $\cA$}
$$
and with these choices we have :
$$
\iota_{f^o}=(\pi_f)^o
\qquad
\pi_{f^o}=(\iota_f)^o
\qquad
\beta_{f^o}=(\beta_f)^o
$$
so $\beta_{f^o}$ is an isomorphism whenever the same holds
for $\beta_f$.

(ii)\ \
Notice that for every category $\cC$, the category
$\bFun(\cC,\cA)$ is pre-additive; indeed, if
$\tau,\sigma:F\Rightarrow G$ are two natural
transformations between functors $F,G:\cC\to\cA$,
then we obtain a natural transformation $\tau+\sigma$ from
$F$ to $G$, by the rule : $(\tau+\sigma)_X:=\tau_X+\sigma_X$
for every $X\in\Ob(\cC)$ (where the sum denotes the addition
law of $\Hom_\cA(FX,GX)$). Clearly, this rule yields an
abelian group structure, and the composition of natural
transformation defines a bilinear pairing
$(\tau,\tau')\mapsto\tau\odot\tau'$ on the resulting groups of
natural transformations (verification left to the reader).

Moreover, if $\cA$ is an additive (resp. abelian) category,
then the same holds for $\bFun(\cC,\cA)$, for every
category $\cC$. Indeed, if $0$ denotes a zero object of
$\cA$, it follows easily from remark \ref{rem_prolyx}(iii)
that the constant functor $0_\cC:\cC\to\cA$ with
$0_\cC(C):=0$ for every $C\in\Ob(\cC)$ is a zero object for
$\bFun(\cC,\cA)$. By the same token, if $F,G:\cC\to\cA$ are
any two functors, then the biproduct $F\oplus G$ is represented
by the functor given by the rule : $C\mapsto FC\oplus GC$
for every $C\in\Ob(\cC)$ (where $FC\oplus GC$ denotes any
fixed choice of an object of $\cA$ representing the biproduct
of $FC$ and $GC$); then $(F\oplus G)(\phi)$ is the induced
morphism $F\phi\oplus G\phi:FC\oplus GC\to FC'\oplus GC'$
for every morphism $\phi:C\to C'$ in $\cC$ (details left
to the reader). Likewise, for every natural transformation
$\tau:F\Rightarrow G$, the kernel and cokernel of $\tau$
are computed argumentwise, and the same holds for the
natural morphism $\beta_\tau$ of \eqref{eq_ker-coker}, so
if all kernels and cokernels are representable in $\cA$,
the same holds for the kernels and cokernels in
$\bFun(\cC,\cA)$, and if $\beta_{\tau_C}$ is an isomorphism
for every $C\in\Ob(\cC)$, then the same holds for
$\beta_\tau$, whence the contention.

(iii)\ \
For any other pre-additive category $\cB$, let us denote
by $\mathbf{Add}(\cB,\cA)$ the full subcategory of
$\bFun(\cB,\cA)$ whose objects are the additive functors.
It follows from (ii) that $\mathbf{Add}(\cB,\cA)$ is a
pre-additive category.

(iv)\ \
Suppose that $\cA$ has small $\Hom$-sets. By definition,
for every $A,B\in\Ob(\cA)$, the set $h_A(B):=\Hom_\cA(B,A)$
carries an abelian group structure, such that the presheaf
$h_A$ factors through an additive functor
$h^\dagger_A:\cA^o\to\Z\Mod$ from $\cA^o$ to the category of
(small) abelian groups, and the forgetful functor
$\Z\Mod\to\Set$. Hence, the Yoneda embedding factors through
a fully faithful {\em group-valued Yoneda embedding}
$$
h^\dagger:\cA\to\mathbf{Add}(\cA^o,\Z\Mod).
$$
In view of (iii), we conclude that every pre-additive category
is a full subcategory of an abelian category. Moreover, Yoneda's
lemma extends {\em verbatim} to the group-valued case : namely,
by inspecting the proof of proposition \ref{prop_yoneda}, we
see that, for every $A\in\Ob(\cA)$ and every {\em additive} functor
$F:\cA^o\to\Z\Mod$ there are natural isomorphisms of abelian groups
\set\begin{equation}\label{eq_Yoneda-additive}
FA\isom\Hom_{\mathbf{Add}(\cA^o,\Z\Mod)}(h^\dagger_A,F).
\end{equation}

(v)\ \
Let $f:\cA\to\cB$ be any functor between pre-additive
categories; then the arguments of remark \ref{rem_was-cofinal}(i)
extend {\em verbatim} to the present situation : after
choosing a universe $\sU$ large enough so that $\cA$ and
$\cB$ are $\sU$-small, $f$ induces functor
$$
f^*:\bFun(\cB^o,\Z\Mod)\to\bFun(\cA^o,\Z\Mod)
$$
that admits both left and right adjoints, denoted respectively
$f_!$ and $f_*$, and we have :
\end{remark}

\begin{proposition} In the situation of remark
{\em\ref{rem_Add.Fun}(v)}, suppose that $f$ is additive
and both $\cA$ and $\cB$ are additive categories. Then :
\begin{enumerate}
\item
Both $f^*$, $f_!$ and $f_*$ are additive functors, and
restrict to functors
$$
\xymatrix{
\mathbf{Add}(\cB^o,\Z\Mod)\ar@<.5ex>[rr]^-{f^*} & &
\mathbf{Add}(\cA^o,\Z\Mod). \ar@<.5ex>[ll]^-{f_!\ \ f_*} 
}$$
\item
The resulting diagram of functors
$$
\xymatrix{
\cA \ar[r]^-{h^\dagger} \ar[d]_f &
\mathbf{Add}(\cA^o,\Z\Mod) \ar[d]^{f_!} \\
\cB \ar[r]^-{h^\dagger} & \mathbf{Add}(\cB^o,\Z\Mod)
}$$
is essentially commutative.
\end{enumerate}
\end{proposition}
\begin{proof}(i): Since every left (resp. right) adjoint
functor is right (resp. left) exact, remark
\ref{rem_fun-exact-add}(iii) says that $f^*$, $f_*$ and
$f_!$ are additive. Next, a simple inspection shows that
$f^*$ transforms additive functors into additive functors.
Let now $F:\cA^o\to\Z\Mod$ be an additive functor,
$B\in\Ob(\cB)$ any object, and $G:=f_!F$; from the
proof of theorem \ref{th_Kan-ext}, we see that
$$
GB=\colim_{\psi:B\to fA}FA
$$
where the colimit ranges over the small category $f\cA^o/B$
of all pairs $(A,\psi)$ consisting of an object $A$ of
$\cA$, and a morphism $\psi:B\to fA$ in $\cB$. Denote by
$0_\cA$ and $0_\cB$ the zero objects of $\cA$ and $\cB$; we
wish to show that $G$ is additive, and according to remark
\ref{rem_fun-exact-add}(iii), it suffices to check that
$G(0_\cB)=0$, and that the natural morphism
$G(B\oplus B)\to GB\oplus GB$ (deduced from the projections
$p_i:B\oplus B\to B$) is an isomorphism, for every $B\in\Ob(\cB)$.

However, notice that the functor $\ss_{0_\cB}:f\cA^o/0_\cB\to\cA^o$
is an isomorphism of categories (notation of
\eqref{subsec_fibreovercat}); whence a natural isomorphism
$$
G(0_\cB)\isom\colim_{\cA^o}F\isom F(0_\cA)=0
$$
where the last identity holds, since $F$ is additive.
Next, for any $B_1,B_2\in\Ob(\cB)$ consider the functor
$$
\Phi:(f\cA^o/B_1)\times(f\cA^o/B_2)\to f\cA^o/B_1\oplus B_2
\qquad
((A_1,\psi_1),(A_2,\psi_2))\mapsto
(A_1\oplus A_2,\psi_1\oplus\psi_2).
$$

\begin{claim}\label{cl_cofinal-cat}
The functor $\Phi$ is cofinal.
\end{claim}
\begin{pfclaim} We apply the criterion of proposition
\ref{prop_MacL-cofinal}. Indeed, let
$$
i:=(A,\psi:B_1\oplus B_2\to fA)
$$
be any object of $f\cA^o/B_1\oplus B_2$ and
$$
\begin{aligned}
j:=\, & ((A_1,\phi_1:B_1\to fA_1),(A_2,\phi_2:B_2\to fA_2)) \\
j':=\, & ((A'_1,\phi'_1:B_1\to fA'_1),(A'_2,\phi'_2:B_2\to fA'_2))
\end{aligned}
$$
any two objects of $(f\cA^o/B_1)\times(f\cA^o/B_2)$,
and suppose that $\beta:i\to\Phi j$ and $\beta':i\to\Phi j'$
are morphisms in $f\cA^o/B_1\oplus B_2$, given by morphisms
$\beta:A_1\oplus A_2\to A$ and $\beta':A_1'\oplus A'_2\to A$
in $\cA$ such that
$f\beta\circ(\phi_1\oplus\phi_2)=\psi=
f\beta'\circ(\phi'_1\oplus\phi'_2)$. We let
$$
j'':=((A_1\oplus A'_1,\delta_1:B_1\to f(A_1\oplus A'_1)),
(A_2\oplus A'_2,\delta_2:B_2\to f(A_2\oplus A'_2))
$$
where $\delta_1:=(\phi_1\oplus\phi'_1)\circ\Delta_{B_1}$,
and likewise for $\delta_2$. Let also
$$
A_1\oplus A_2\xleftarrow{\ p\ }
A_1\oplus A'_1\oplus A_2\oplus A'_2
\xrightarrow{\ p'\ }A'_1\oplus A'_2
$$
be the natural projections, which define morphisms
$$
\Phi p:\Phi j\to\Phi j''
\qquad
\Phi p':\Phi j'\to\Phi j''
\qquad
\text{in $f\cA^o/B_1\oplus B_2$}
$$
A simple inspection shows that
$\Phi(p)\circ\beta=\Phi(p')\circ\beta'$, and therefore
$i/\Phi(f\cA^o/B_1\oplus B_2)$ is connected, whence
the claim.
\end{pfclaim}

In light of claim \ref{cl_cofinal-cat}, we are reduced to
checking that the natural morphism
$$
\colim_{(f\cA^o/B_1)\times(f\cA^o/B_2)}
F\circ\ss_{B_1\oplus B_2}\circ\Phi\to GB_1\oplus GB_2
$$
is an isomorphism, for any $B_1,B_2\in\Ob(\cB)$. The latter
assertion follows easily by inspecting the definitions,
since $F$ is additive. Lastly, a similar argument shows
that $f_*F$ is additive, whenever the same holds for $F$ :
the reader can spell out the proof as an exercise.

(ii): One may argue as in remark \ref{rem_was-cofinal}(ii) :
in view of \eqref{eq_Yoneda-additive}, we see that, for
every object $A$ of $\cA$, both $h^\dagger_{fA}$ and
$f_!h^\dagger_A$ represent the same functor : details
left to the reader.
\end{proof}

\begin{lemma}\label{lem_abel-fincomplete}
For any abelian category $\cA$, the following holds :
\begin{enumerate}
\item
$\cA$ is finitely complete and finitely cocomplete.
\item
The image of any morphism $f:A\to B$ of $\cA$ exists
in $\cA$, and we have
$$
\Coker\,\pi_f=\Img\,f
\qquad
\text{in\ \ $\mathrm{Sub}(B)$}
$$
(notation of \eqref{subsec_well-power}, example
{\em\ref{ex_pre-misc}} and remark {\em\ref{rem_additive-cat}(v)}).
\end{enumerate} 
\end{lemma}
\begin{proof}(i): By proposition
\ref{prop_complete-criteria}(i,ii), it suffices to
check that all equalizers and coequalizers are
representable in $\cA$. However, let $g,g':X\to Y$
be any pair of morphisms in $\cA$ with same soruces
and targets; it easily seen that $\Ker\,(g-g')$
(resp. $\Coker\,(g-g')$) represents the equalizer
(resp. the coequalizer) of $g$ and $g'$, whence
the assertion.

(ii): Suppose that $f$ factors through a subobject
$g:C\to B$ of $B$; there follows a natural morphism
$$
\Coker\,f\to\Coker\,g
$$
which in turns induces a morphism
$\Ker\,\pi_f\to\Ker\,\pi_g$ of subobjects of $B$.
But since $\cA$ is abelian, $\Ker\,\pi_g$ is also
the cokernel of $\iota_g:\Ker\,g\to C$ (notation of
remark \ref{rem_additive-cat}(v)). However, $\Ker\,g=0$
by remark \ref{rem_additive-cat}(vi), hence
$\Coker\,\iota_g=C$, whence the contention.
\end{proof}

\begin{definition} Let $\cA$ be any abelian category
with zero object $0\in\Ob(\cA)$, and
$$
0\xrightarrow{\ 0\ }A\xrightarrow{\ f\ }
B\xrightarrow{\ g\ } C\xrightarrow{\ 0\ } 0
$$
a sequence of morphisms in $\cA$.
\begin{enumerate}
\item
We say that the sequence $(f,g)$ is {\em exact} if
$\Ker\,g=\Img\,f$.
\item
We say that the sequence $(0,f,g)$ is {\em left
exact} if both $(0,f)$ and $(f,g)$ are exact.
\item
We say that the sequence $(f,g,0)$ is {\em right
exact} if both $(f,g)$ and $(g,0)$ are exact.
\item
We say that the sequence $(0,f,g,0)$ is {\em short
exact} if $(0,f,g)$ is left exact and $(f,g,0)$ is
right exact.
\end{enumerate}
\end{definition}
In other words, the sequence $(0,f,g)$ is left exact
if and only if $f$ is a monomorphism, and $A$ represents
$\Ker\,g$. Likewise, $(f,g,0)$ is right exact if and only
if $g$ is an epimorphism and the image of $f$ represents
$\Ker\,g$. The terminology is explained by the following :

\begin{lemma} Let $F:\cA\to\cB$ be any additive
functor between abelian categories. Then :
\begin{enumerate}
\item
$F$ is left exact if and only if it transforms
left exact sequences of $\cA$ into left exact
sequences of $\cB$.
\item
$F$ is right exact if and only if it transforms
right exact sequences of $\cA$ into right exact
sequences of $\cB$.
\item
$F$ is exact if and only if it transforms short
exact sequences of $\cA$ into short exact sequences
of $\cB$.
\end{enumerate}
\end{lemma}
\begin{proof}(i): Any left exact functor transforms
kernels into kernels and monomorphisms into
monomorphisms (proposition \ref{prop_was-also-cofinal}(i)),
so the condition is necessary. Conversely, suppose
that $F$ fulfills this condition; since $F$ is
additive, it commutes with finite products, so
it suffices to check that it also commutes with
equalizers (proposition \ref{prop_commute-criteria}(i)).
Arguing as in the proof of lemma \ref{lem_abel-fincomplete}(i),
we reduce to checking that $F$ commutes with kernels.
But this follows by considering the left exact
sequence
$$
0\to\Ker\,f\to A\xrightarrow{\ f\ }B
\qquad
\text{for any morphism $f$ of $\cA$}
$$
and its image under $F$. Assertion (ii) admits the
dual proof, and (iii) follows by considering the
similar short exact sequences
$$
0\to\Ker\,f\to A\to\Img\,f\to 0
\qquad
0\to\Img\,f\to B\to\Coker\,f\to 0
$$
for every $f$ as in the foregoing, and arguing
as in the proof of (i).
\end{proof}

\begin{definition}\label{def_ab-tensor-cat}
An {\em abelian tensor category} is a tensor category
$(\cC,\otimes,\Phi,\Psi)$ such that $\cC$ is an abelian
category, and the functor $\otimes$ induces bilinear pairings
$$
\Hom_\cC(A,B)\times\Hom_\cC(A',B')\to\Hom_\cC(A\otimes A',B\otimes B')
\quad :\quad
(f,g)\mapsto f\otimes g
$$
for every $A,A',B,B'\in\Ob(\cC)$.
\end{definition}

\begin{remark} Let $(\cA,\otimes,\Phi,\Psi)$ be a tensor
category, such that $\cA$ is abelian. If $\cA$ admits an
internal $\Hom$ functor $\cHom$, the functor $-\otimes A$
is right exact, and the functor $\cHom(A,-)$ is left
exact for every $A\in\Ob(\cA)$, so both are additive,
by virtue of remark \ref{rem_fun-exact-add}(iii).
Especially, $\cA$ is an abelian tensor category, in this
case.
\end{remark}

\begin{lemma}\label{lem_simple-imbeddings}
Let $\cA$ be any additive category with small $\Hom$-sets,
and $\Sigma\subset\Ob(\cA)$ a small subset. We have :
\begin{enumerate}
\item
If $\cA$ is abelian, there exists a small full abelian
subcategory $\cB$ of $\cA$ such that $\Sigma\subset\Ob(\cB)$.
\item
If $\cA$ is small, there exist a complete and cocomplete
abelian tensor category $(\cC,\otimes)$ with internal $\Hom$
functor, and a fully faithful additive functor $\cA\to\cC$. 
\end{enumerate}
\end{lemma}
\begin{proof}(i): Let $\cB_0$ be the full subcategory of
$\cA$ such that $\Ob(\cB_0)=\Sigma$; clearly $\cB_0$ is
small. Next, for any subcategory $\cD$ of $\cA$, denote
by $\cD'$ a subcategory of $\cA$ obtained as follows.
For every morphism $\phi$ of $\cD$, we pick objects in
$\cA$ representing the kernel and cokernel of $\phi$,
and for any two objects of $\cD$, we pick an object in
$\cA$ representing their product; let $\Sigma'\subset\Ob(\cA)$
be the resulting subset. Then $\cD'$ is the full subcategory
of $\cA$ such that $\Ob(\cD')=\Ob(\cD)\cup\Sigma'$. It is
easily seen that $\cD'$ is small, whenever the same holds
for $\cD$. Then we set inductively $\cB_{i+1}:=\cB'_i$ for
every $i\in\N$. The full subcategory $\cB$ of $\cA$ with
$\Ob(\cB)=\bigcup_{i\in\N}\Ob(\cB_i)$ is still small, and
it is abelian, by construction.

(ii): We let $\cC:=\bFun(\cA,\Z\Mod)$. Then $\cC$ is an
abelian category, by virtue of remark \ref{rem_Add.Fun}(ii),
and since $\Z\Mod$ is complete and cocomplete, the same holds
for $\cC$; moreover, the standard tensor product of abelian
groups defines a tensor category structure with internal $\Hom$
functor on $\Z\Mod$, and the latter is inherited by $\cC$
(remarks \ref{rem_Mac-Lane}(ii) and \ref{rem_Hom-how-to}(vi)).
It is clear that these two structures amount to an abelian
tensor category, and the group-valued Yoneda embedding is
the sought fully faithful functor.
\end{proof}

\sset\subsubsection{}\label{subsec_mixed-tensors}
Let $\cA$ be an additive category with small $\Hom$-sets,
and
$$
h^\dagger:\cA^o\to\cA^\dagger:=\bFun(\cA,\Z\Mod)
$$
the fully faithful group-valued Yoneda embedding. For
every abelian group $G$, denote by $G_\cA:\cA\to\Z\Mod$
the constant functor with value $G$ : so, $G_\cA(A):=A$
for every $A\in\Ob(\cA)$, and $G_\cA(\phi):=\one_A$ for
every morphism $\phi$ in $\cA$. Since $\cA^\dagger$ is an
abelian tensor category with an internal $\Hom$ functor
$\cHom$ (see the proof of lemma
\ref{lem_simple-imbeddings}(ii)), we may define
$$
G\otimes_\Z A:=\cHom(G_\cA,h^\dagger_{A^o})
\qquad
\text{for every $A\in\Ob(\cA)$}.
$$
If $G$ is free of finite rank, $G\otimes_\Z A$ is a finite
direct sum of copies of $A$, and therefore lies in the
essential image of $h^\dagger$; the same holds for a
finitely generated $G$, provided $\cA$ is an abelian
category : indeed, we may write $G$ as a cokernel of
a map $L_1\to L_2$ of free abelian groups of finite rank,
and since the functor $\cHom(-,h^\dagger_A)$ is right exact,
we see that $G\otimes_\Z A$ is the kernel of a morphism
of $\cA^o$ ({\em i.e.} the cokernel of a morphism of
$\cA$), so it is represented by an object of $\cA$.
Moreover, if $\phi:G\to H$ is any morphism of abelian
groups, we have an obvious induced morphism
$\phi_\cA:G_\cA\to H_\cA$, whence a morphism
$\phi\otimes_\Z A:=\cHom(\phi_\cA,h^\dagger_A)$.

Thus, if $\cA$ is abelian, after replacing
$G\otimes_\Z A$ by an isomorphic object, we
obtain a well defined functor
\set\begin{equation}\label{eq_mixed-tensor}
\Z\Mod_\mathrm{fg}\times\cA\to\cA
\qquad
(G,A)\mapsto G\otimes_\Z A
\end{equation}
where $\Z\Mod_\mathrm{fg}$ is the full subcategory of $\Z\Mod$
whose objects are the finitely generated abelian groups.
This functor is not unique, but any two such functors are
naturally isomorphic. If $\cA$ is only additive, we can
still define such a tensor product functor on the
category $\Z\Mod_\mathrm{fft}\times\cA\to\cA$, where
$\Z\Mod_\mathrm{fft}$ is the full subcategory of
$\Z\Mod_\mathrm{fg}$ whose objects are the free abelian
groups of finite rank.

\begin{remark}\label{rem_mixed-tensors}
Keep the notation of \eqref{subsec_mixed-tensors}; we have :

(i)\ \
From the construction, it is clear that \eqref{eq_mixed-tensor}
is a {\em biadditive} functor, {\em i.e.}, for every abelian group
$G$, and every $A\in\Ob(\cA)$, the restrictions $G\otimes-$ and
$-\otimes A$ of \eqref{eq_mixed-tensor} are additive.

(ii)\ \
Suppose that $\cA$ is cocomplete; since the tensor product
is right exact, it follows easily that \eqref{eq_mixed-tensor}
extends to the whole of $\Z\Mod$ : details left to the reader.
\end{remark}

\section{Sites and topoi}
In this chapter, we assemble some generalities concerning sites
and topoi. The main reference for this material is \cite{SGA4-1}.

\subsection{Topologies and sites}\label{sec_topoi}
As in the previous chapter, we fix a universe $\sU$ such that
$\N\in\sU$, and small means $\sU$-small throughout. Especially,
a presheaf on any category takes its values in $\sU$, unless
explicitly stated otherwise.

\begin{definition}\label{def_topology}
Let $\cC$ be a category.
\begin{enumerate}
\item
A {\em topology\/} on $\cC$ is the datum, for every $X\in\Ob(\cC)$,
of a set $J(X)$ of sieves of $\cC\!/\!X$, fulfilling the following
conditions :
\begin{enumerate}
\item
(Stability under base change)\ \ For every morphism $f:Y\to X$ of
$\cC$, and every $\cS\in J(X)$, the sieve $\cS\times_Xf$ lies in
$J(Y)$.
\item
(Local character)\ \ Let $X$ be any object of $\cC$, and $\cS$,
$\cS'$ two sieves of $\cC\!/\!X$, with $\cS\in J(X)$. Suppose that,
for every object $f:Y\to X$ of $\cS$, the sieve $\cS'\times_Xf$ lies
in $J(Y)$. Then $\cS'\in J(X)$.
\item
For every $X\in\Ob(\cC)$, we have $\cC\!/\!X\in J(X)$.
\end{enumerate}
\item
In the situation of (i), the elements of $J(X)$ shall be called the
{\em sieves covering $X$}. Moreover, say that $\cS$ is the sieve of
$\cC\!/\!X$ generated by a family $(f_i:X_i\to X~|~i\in I)$ of
morphisms. If $\cS$ is a sieve covering $X$, we say that the family
$(f_i~|~i\in I)$ {\em covers $X$\/}, or that it is a {\em covering
family of $X$}.
\item
The datum $(\cC,J)$ of a category $\cC$ and a topology
$J:=(J(X)~|~X\in\Ob(\cC))$ on $\cC$ is called a {\em site\/}, and
then $\cC$ is also called the {\em category underlying\/} the site
$(\cC,J)$. We say that $(\cC,J)$ is a {\em small site}, if $\cC$
is a small category.
\item
The set of all topologies on $\cC$ is partially ordered by
inclusion: given two topologies $J_1$ and $J_2$ on $\cC$,
we say that $J_1$ is {\em finer\/} than $J_2$, if
$J_2(X)\subset J_1(X)$ for every $X\in\Ob(\cC)$.
\end{enumerate}
\end{definition}

\begin{example}\label{ex_top-spaces}
Let $T$ be any topological space, and $\cT$ the set of
open subsets of $T$, which is partially ordered by
inclusion, and can thus be regarded naturally as a
category, as in example \ref{ex_universe}(iii).
Then the category $\cT$ carries a natural topology $J_T$
in the sense of definition \ref{def_topology} : namely,
for every $U\in\cT$ we declare that the elements of
$J_T(U)$ are the sieves $\cS\subset\cT/U$ such that
$\bigcup_{(f:U'\to U)\in\cS}U'=U$. In other words, a family
$(U'_i\to U~|~i\in I)$ of morphisms of $\cT$ covers $U$
for the topology $J_T$ if and only if $\bigcup_{i\in I}U_i=U$.
\end{example}

\begin{remark}\label{rem_topology}
Let $(\cC,J)$ be any site.

(i)\ \
Any finite intersection of sieves covering an object $X$, again
covers $X$. Indeed, say that $\cS_1$ and $\cS_2$ are two sieves
covering $X$; set $\cS:=\cS_1\cap\cS_2$ and let $f:Y\to X$ be any
object of $\cS_1$. Then $\cS\times_Xf=\cS_2\times_Xf\in J(Y)$, so
the assertion follows from the local character of $J$.

(ii)\ \
Also, any sieve of $\cC\!/\!X$ containing a covering sieve
is again a covering sieve. Indeed, if $\cS\subset\cS'$, then
$\cS'\times_Xf=\cC/Y$ for every object $f:Y\to X$ of $\cS$.

(iii)\ \
Let $f_\bullet:=(f_i:Y_i\to X~|~i\in I)$ be a family of objects
of $\cC\!/\!X$ that generates a sieve $\cS$ covering $X$, and
for every $i\in I$, let $(g_{ij}:Z_{ij}\to Y_i~|~i\in J_i)$ be
a family of objects of $\cC/Y_i$ that generates a sieve $\cS_i$
covering $Y_i$. Then the family
$(f_i\circ g_{ij}:Z_{ij}\to X~|~i\in I,\ j\in J_i)$ is a
refinement of $f_\bullet$ (see definition \ref{def_sieve}(ii))
and generates a sieve $\cS'$ covering $X$. Indeed, say that
$f:Y\to X$ lies in $\cS$, and pick $i\in I$ such that $f$
factors through $f_i$ and a morphism $g:Y\to Y_i$; then it
is easily seen that $\cS_i\times_{Y_i}g\subset\cS'\times_Xf$.
\end{remark}

\sset\subsubsection{}\label{subsec_alt-topolo}
Suppose that $\cC$ has small $\Hom$-sets. Then, in view
of remark \ref{rem_sieves-and-sub}(ii), a topology can also
be defined by assigning, to any object $X$ of $\cC$, a family
$J'(X)$ of subobjects of $h_X$, called the
{\em subobjects covering\/} $X$, such that :
\begin{enumerate}
\alphaenu
\item
For every $X\in\Ob(\cC)$, every $R\in J'(X)$, and every morphism
$Y\to X$ in $\cC$, the fibre product $R\times_XY$ lies in $J'(Y)$.
\item
Say that $X\in\Ob(\cC)$, and let $R$, $R'$ be two subobjects of
$h_X$, such that $R\in J'(X)$. Suppose that, for every
$Y\in\Ob(\cC)$, and every morphism $f:h_Y\to R$, we have
$R'\times_XY\in J'(Y)$. Then $R'\in J'(X)$.
\item
$h_X\in J'(X)$ for every $X\in\Ob(\cC)$.
\end{enumerate}
This viewpoint is adopted in the following :

\begin{definition}\label{def_sheaf}
Let $\sV$ be a universe, $C:=(\cC,J)$ a site, and
$F\in\Ob(\cC^\wedge_\sV)$.
\begin{enumerate}
\item
We say that $F$ is a {\em separated $\sV$-presheaf} (resp. a
{\em $\sV$-sheaf}) on $C$, if for every $X\in\Ob(\cC)$, every
subobject $R$ covering $X$, and every universe $\sV'$ containing
$\sV$ and such that the category $\cC$ has $\sV'$-small $\Hom$-sets,
the induced morphism :
$$
F(X)=\Hom_{\cC^\wedge_{\sV'}}(h_X,F)\to\Hom_{\cC^\wedge_{\sV'}}(R,F)
$$
is injective (resp. is bijective). For $\sV=\sU$, we shall
just say separated presheaf instead of separated $\sU$-presheaf,
and likewise for sheaves.
\item
We let $C^\sim_\sV$ (resp. $C^\sep_\sV$) be the full subcategory
of $\cC^\wedge_\sV$ whose objects are the $\sV$-sheaves (resp. the
separated $\sV$-presheaves) on $C$. For $\sV=\sU$, this category
will be usually denoted simply $C^\sim$ (resp. $C^\sep$).
\end{enumerate}
\end{definition}

\sset\subsubsection{}\label{subsec_interpret-descent}
In the situation of definition \ref{def_sheaf}(i), say that
$R=h_\cS$ for some sieve $\cS$ covering $X$. In light of
\eqref{eq_colim-sieve}, \eqref{eq_colim-and-Hom} and proposition
\ref{prop_yoneda}(ii), we get a natural identification :
$$
\Hom_{\cC^\wedge_{\sV'}}(R,F)\isom
\Hom_{\cC^\wedge_{\sV'}}(\colim_\cS h_\cC\circ\ss_\cS,F)\isom
\lim_{\cS^o} F\circ\ss^o_\cS.
$$
Let also $(S_i\to X~|~i\in I)$ be a generating family for $\cS$.
Combining lemma \ref{lem_coeq-sieve}, example
\ref{ex_fibred-cats-II}(i) and \eqref{subsec_fibred-cats-II},
we get a natural isomorphism :
$$
\Hom_{\cC^\wedge_{\sV'}}(R,F)\isom
\Equal\Bigl( \prod_{i\in I}
\xymatrix{\sCart_\cC(\cC\!/\!S_i,\cFib(F))
\ar@<-.5ex>[r] \ar@<.5ex>[r] &} \!\!\!\!\!\!\!
\prod_{(i,j)\in I\times I}\!\!\!\sCart_\cC(\cC\!/\!S_{ij},\cFib(F))
\Bigr)
$$
which -- again by \eqref{subsec_fibred-cats-II} --
we may rewrite more simply as :
$$
\Hom_{\cC^\wedge_{\sV'}}(R,F)\isom
\Equal\Bigl( \prod_{i\in I}
\xymatrix{F(S_i)\ar@<-.5ex>[r] \ar@<.5ex>[r] &} \!\!\!\!\!\!\!
\prod_{(i,j)\in I\times I}\!\!\!
\Hom_{\cC^\wedge_{\sV'}}(h_{S_i}\times_{h_X}h_{S_j},F)
\Bigr).
$$
The above equalizer can be described explicitly as follows.
It consists of all the systems
$$
(a_i~|~i\in I)
\qquad
\text{with $a_i\in F(S_i)$\ \ for every $i\in I$}
$$
such that, for every $i,j\in I$ and every object $Y\to X$ of
$\cC\!/\!X$, we have :
\set\begin{equation}\label{eq_interpet-descent}
F(g_i)(a_i)=F(g_j)(a_j)
\qquad
\text{for every pair $(S_i\xleftarrow{g_i}Y\xrightarrow{g_j}S_j)$
of morphisms in $\cC\!/\!X$}.
\end{equation}
If all $S_{ij}:=S_i\times_XS_j$ are representable in $\cC$, the latter
expression takes the more familiar form:
$$
\Hom_{\cC^\wedge_{\sV'}}(R,F)\isom
\Equal\Bigl( \prod_{i\in I}
\xymatrix{F(S_i)\ar@<-.5ex>[r] \ar@<.5ex>[r] &} \!\!\!\!\!\!\!
\prod_{(i,j)\in I\times I}\!\!\!F(S_{ij})
\Bigr).
$$

\begin{remark}\label{rem_sheaves}
Let $C:=(\cC,J)$ be a site such that $\cC$ is a category with small
$\Hom$-sets.

(i)\ \ 
The arguments in \eqref{subsec_interpret-descent} yield also the
following. A presheaf $F$ on $\cC$ is separated (resp. is a sheaf)
on $C$, if and only if every covering sieve of $C$ is a sieve of
$1$-descent (resp. of $2$-descent) for the fibration
$\ss_F:\cFib(F)\to\cC$ of \eqref{subsec_fibred-cats-II}.

(ii)\ \
Let $F:\cA\to\cB$ be a fibration between two categories,
and $i\leq 2$ an integer.
For every $X\in\Ob(\cB)$, let $J^i_F(X)$ denote the set of
all sieves $\cS\subset\cB/X$ of universal $F$-$i$-descent.
Then we claim that $J^i_F$ is a topology on $\cB$. Indeed, it is
clear that $J^i_F$ fulfills conditions (a) and (c) of definition
\ref{def_topology}(i), and condition (b) follows from proposition
\ref{prop_topol-of-univ-descent}.

(iii)\ \ 
Let $F$ be a presheaf on $\cC$. We deduce from (i) and (ii)
that the topology $J^F:=J^2_{\ss_F}$ is the finest on $\cC$
for which $F$ is a sheaf. A subobject $R\subset h_X$ (for
any $X\in\Ob(\cC)$) lies in $J^F(X)$ if and only if the
natural map $F(X')\to\Hom_{\cC^\wedge}(R\times_XX',F)$ is
bijective for every morphism $X'\to X$ in $\cC$.

(iv)\ \ 
More generally, if $(F_i~|~i\in I)$ is any family of
presheaves on $\cC$, there exists a finest topology for
which each $F_i$ is a sheaf : namely, the intersection
of the topologies $J^{F_i}$ as in (iii).

(v)\ \ 
As an important special case, we deduce the existence
of a finest topology $J$ on $\cC$ such that all representable
presheaves are sheaves on $(\cC,J)$. This topology $\Can_\cC$
is called the {\em canonical topology} on $\cC$. We thus
associate to every category $\cC$ a {\em canonical site}
$$
\Can(\cC):=(\cC,\Can_\cC).
$$

(vi)\ \
Another interesting case of (iv) is obtained by taking
the family of all presheaves on $\cC$. The corresponding
topology $J$ on $\cC$ can be easily described explicitly :
namely, one takes $J(X):=\{\cC/X\}$ for every $X\in\Ob(\cC)$.
\end{remark}

\begin{example}\label{ex_strict-epis}
Let $\cC$ be a category, $X\in\Ob(\cC)$ any object, and
$\cS\subset\cC/X$ a sieve.

(i)\ \
We say that $\cS$ is an {\em epimorphic sieve} (resp. a
{\em strict epimorphic sieve}) if for every $Y\in\Ob(\cC)$
the natural map
$$
h_Y(X)\to\Hom_{\cC^\wedge}(h_\cS,h_Y)
$$
is injective (resp. bijective). We say that $\cS$ is a
{\em universal epimorphic sieve} (resp. a {\em universal
strict epimorphic sieve}) if $\cS\times_Xh$ is epimorphic
(resp. strict epimorphic) for every morphism $h:Y\to X$ of
$\cC$. From remark \ref{rem_sheaves}(iii,v) we see that $\cS$
is universal strict epimorphic if and only if it covers $X$
in the canonical topology of $\cC$.

(ii)\ \
We say that a family of morphisms
$f_\bullet:=(f_i:X_i\to X~|~i\in I)$ of $\cC$ is {\em epimorphic}
(resp. {\em strict epimorphic}) if it generates an epimorphic
(resp. strictly epimorphic) sieve. This is the same as saying
for every $Y\in\Ob(\cC)$, the natural map
\set\begin{equation}\label{eq_strict-epimorph}
\Hom_\cC(X,Y)\to\prod_{i\in I}\Hom_\cC(X_i,Y)
\qquad
g\mapsto(g\circ f_i~|~i\in I)
\end{equation}
is injective (resp. and its image consists of all the systems
of morphisms $(g_i:X_i\to Y~|~i\in I)$ such that, for every
$Z\in\Ob(\cC)$, every $i,j\in I$, and every pair of morphisms
$X_j\xleftarrow{h_j}Z\xrightarrow{h_i}X_i$ with
$f_i\circ h_i=f_j\circ h_j$, we have $g_i\circ h_i=g_j\circ h_j$).

(iii)\ \
We say that a family $f_\bullet$ as in (ii) is {\em effective
epimorphic}, if it is strict epimorphic, and moreover the fibre
products $X_i\times_XX_j$ are representable in $\cC$ for every
$i,j\in I$. This is the same as saying that for every $Y\in\Ob(\cC)$
the map \eqref{eq_strict-epimorph} identifies $\Hom_\cC(X,Y)$ with
the equalizer of the two natural maps
$$
\prod_{i\in I}\xymatrix{\Hom_\cC(X_i,Y) \ar@<.5ex>[r] \ar@<-.5ex>[r]
& }\!\!\!\!\!\prod_{(i,j)\in I\times I}\!\!\!\!\Hom_\cC(X_i\times_XX_j,Y).
$$

(iv)\ \
A family $f_\bullet$ as in (ii) is called {\em universal epimorphic}
(resp. {\em universal effective epimorphic}) if (a) the fibre products
$Y_i:=X_i\times_XY$ are representable in $\cC$, for every $i\in I$
and every morphism $Y\to X$ in $\cC$, and (b) all the resulting
families $(Y_i\to Y~|~i\in I)$ are still epimorphic (resp.
effective epimorphic). Then clearly such a family generates
a universal epimorphic sieve (resp. a universal strict epimorphic
sieve).

(v)\ \
We say that a morphism $f:X'\to X$ of $\cC$ is an
{\em effective epimorphism} (resp. a {\em universal epimorphism},
resp. a {\em universal effective epimorphism}) if the family
$\{f\}$ has the corresponding property. Hence, $f$ is an effective
epimorphism if the fibre product $X'\times_XX'$ is representable
in $\cC$, and the induced diagram
$$
\xymatrix{X'\times_XX' \ar@<.5ex>[r]^-{p_1} \ar@<-.5ex>[r]_-{p_2}
& X' \ar[r]^-{f} & X}
$$
identifies $X$ with the coequalizer of the two projections
$p_1$ and $p_2$. (Notice that this condition implies that
$f$ is an epimorphism).
Likewise, $f$ is a universal epimorphism (resp. a universal
effective epimorphism) if for every morphism $Y\to X$,
the fibre product $X'\times_XY$ is representable in $\cC$
and $f\times_XY:X'\times_XY\to Y$ is an epimorphism
(resp. an effective epimorphism).

(vi)\ \
Suppose that all fibre products are representable in $\cC$;
in this case, proposition \ref{prop_diagonal-monos}(ii)
implies that a morphism $X'\to X$ of\/ $\cC$ is a universal
epimorphism if and only if it is an epimorphism and the
coproduct $X\amalg_{X'}X$ is a universal colimit, in the
sense of example \ref{ex_universal-col}(ii).

(vii)\ \
Combining (vi) with example \ref{ex_universal-col}(iv),
we conclude that all epimorphisms of $\Set$ are universal.
Also, it is easily seen that all epimorphisms of $\Set$
are effective, and hence universal effective.
More generally, If\/ $\cC$ is a small category, example
\ref{ex_big-complete}(i) and corollary \ref{cor_pre-misc}(ii)
imply that all epimorphisms in $\cC^\wedge$ are universal effective.
\end{example}

\sset\subsubsection{}\label{subsec_asso-topoi}
Let $C:=(\cC,J)$ be a small site; directly from definition
\ref{def_sheaf} (and from \eqref{eq_lim-and-Hom}), we see
that the category $C^\sim$ is complete, and the fully faithful
inclusion $C^\sim\to\cC^\wedge$ commutes with all limits. Moreover,
given a presheaf $F$ on $\cC$, it is possible to construct a
solution set for $F$ relative to this functor, and therefore
one may apply theorem \ref{th_adjoint-fctr-th} to produce a
left adjoint.
However, a more direct and explicit construction of the
left adjoint can be given; the latter also provides some
additional information which is hard to extract from the
former method. Namely, for every $X\in\Ob(\cC)$, endow $J(X)$
with the partial ordering induced by inclusion of sieves.
Notice that $J(X)$ is small and cofiltered, for every such $X$,
hence the opposite ordered set $J(X)^o$ is filtered.
For a given presheaf $F$ on $\cC$, set 
$$
F^+(X):=\colim_{\cS\in J(X)^o}\Hom_{\cC^\wedge}(h_\cS,F)
$$
where the transition map
$\Hom_{\cC^\wedge}(h_\cS,F)\to\Hom_{\cC^\wedge}(h_{\cS'},F)$ is induced
by the monomorphism $h_{\cS'}\to h_\cS$ of subobjects of $h_X$, for
every inclusion $\cS'\subset\cS$ of sieves covering $X$. Let
$(\tau_\cS:\Hom_{\cC^\wedge}(h_\cS,F)\to F^+(X)\ |\ \cS\in J(X))$
be the universal cocone. For a morphism $f:Y\to X$
in $\cC$ and a sieve $\cS\in J(X)$, the natural projection
$h_\cS\times_XY\to h_\cS$ induces a map
$$
\Hom_{\cC^\wedge}(h_\cS,F)\to\Hom_{\cC^\wedge}(h_\cS\times_XY,F)
\xrightarrow{\ \tau_{\cS\times_Xf}\ }F^+(Y)
$$
whence a map $F^+f:F^+(X)\to F^+(Y)$, after taking colimits.
Likewise, every morphism $F\to F'$ of presheaves on $\cC$
induces a map $\Hom_{\cC^\wedge}(h_\cS,F)\to\Hom_{\cC^\wedge}(h_\cS,F')$,
for every $X\in\Ob(\cC)$ and every $\cS\in J(X)$, whence a map
$F^+(X)\to F'^+(X)$, after taking colimits. Thus, we have a functor :
\set\begin{equation}\label{eq_plus-construction}
\cC^\wedge\to\cC^\wedge
\qquad
F\mapsto F^+.
\end{equation}
with a natural transformation $F(X)\to F^+(X)$, since
$\cC\!/\!X\in J(X)$ for every $X\in\Ob(\cC)$.

\begin{theorem}\label{th_ass-sheaf}
In the situation of \eqref{subsec_asso-topoi}, the following
holds :
\begin{enumerate}
\item
The inclusion functor $i:C^\sim\to\cC^\wedge$ admits the
left adjoint
\set\begin{equation}\label{eq_asso-funct}
\cC^\wedge\to C^\sim
\qquad F\mapsto F^a:=(F^+)^+.
\end{equation}
For every presheaf $F$, we call $F^a$ the
{\em sheaf associated with $F$}.
\item
Morever, the functor \eqref{eq_asso-funct} is exact.
\end{enumerate}
\end{theorem}
\begin{proof} We begin with the following :

\begin{claim}\label{cl_better-injectivity}
Let $F$ be a separated presheaf on $\cC$, and $\cS_1\subset\cS_2$
two sieves covering some $X\in\Ob(\cC)$. Then the natural map
$\Hom_{\cC^\wedge}(h_{\cS_2},F)\to\Hom_{\cC^\wedge}(h_{\cS_1},F)$
is injective.
\end{claim}
\begin{pfclaim} We may find a generating family
$(f_i:S_i\to X~|~i\in I_2)$ for $\cS_2$, and a subset
$I_1\subset I_2$, such that $(f_i~|~i\in I_1)$ generates
$\cS_1$. Let $s,s':h_{\cS_2}\to F$, whose images agree in
$\Hom_{\cC^\wedge}(h_{\cS_1},F)$.
By \eqref{subsec_interpret-descent}, $s$ and $s'$ correspond
to families $(s_i~|~i\in I_2)$, $(s'_i~|~i\in I_2)$ with
$s_i,s'_i\in F(S_i)$ for every $i\in I_2$, fulfilling
the system of identities \eqref{eq_interpet-descent},
and the foregoing condition means that $s_i=s'_i$ for
every $i\in I_1$. We need to show that $s_i=s'_i$ for
every $i\in I_2$. Hence, let $i\in I_2$ be any element;
by assumption, the natural map
$$
F(S_i)\to\Hom_{\cC^\wedge}(h_{\cS_1}\times_XS_i,F)
$$
is injective. However, the objects of $\cS_1\times_Xf_i$
are all the morphisms $g_i:Y\to S_i$ in $\cC$ such that
$f_i\circ g_i=f_j\circ g_j$ for some $j\in I_1$ and some
$g_j:Y\to S_j$ in $\cC$. If we apply the identities
\eqref{eq_interpet-descent} to these maps $g_i$, $g_j$,
we deduce that :
$$
F(g_i)(s_i)=F(g_j)(s_j)=F(g_j)(s'_j)=F(g_i)(s'_i).
$$
In other words, $s_i$ and $s'_i$ have the same image in
$\Hom_{\cC^\wedge}(h_{\cS_1}\times_XS_i,F)$, hence they
agree, as claimed.
\end{pfclaim}

\begin{claim}\label{cl_plus-construction}
(i) The functor \eqref{eq_plus-construction} is left exact.
\begin{enumerate}
\addenu
\item
For every $F\in\Ob(\cC^\wedge)$, the presheaf $F^+$ is separated.
\item
If $F$ is a separated presheaf on $\cC$, then $F^+$ is
a sheaf on $C$.
\end{enumerate}
\end{claim}
\begin{pfclaim} (i) is clear, since $J(X)^o$ is filtered
for every $X\in\Ob(\cC)$.

(ii): Let $s,s'\in F^+(X)$, and suppose that the images of
$s$ and $s'$ agree in $\Hom_{\cC^\wedge}(h_\cS,F^+)$, for
some $\cS\in J(X)$. We may find a sieve $\cT\in J(X)$ such
that $s$ and $s'$ come from elements
$\bar s,\bar s{}'\in\Hom_{\cC^\wedge}(h_\cT,F)$.
Let $(g_i:S_i\to X~|~i\in I)$ be a family of generators
for $\cS$; in view of \eqref{subsec_interpret-descent}, the
images of $s$ and $s'$ agree in $F^+(S_i)$ for every $i\in I$.
The latter means that, for every $i\in I$, there exists
$\cS_i\in J(S_i)$, refining $\cT\times_Xg_i$, such that
the images of $\bar s$ and $\bar s{}'$ agree in
$\Hom_{\cC^\wedge}(h_{\cS_i},F)$. For every $i\in I$,
let $(g_{i\lambda}:T_{i\lambda}\to S_i~|~\lambda\in\Lambda_i)$
be a family of generators for $\cS_i$, and consider the sieve
$\cT'$ of $\cC\!/\!X$ generated by
$(g_i\circ g_{i\lambda}:T_{i\lambda}\to X~|~i\in I,\lambda\in\Lambda_i)$.
Then $\cT'$ covers $X$ (remark \ref{rem_topology}(iii)) and
refines $\cT$, and the images of $\bar s$ and $\bar s{}'$
agree in $\Hom_{\cC^\wedge}(h_{\cT'},F)$ (as one sees easily,
again by virtue of \eqref{subsec_interpret-descent}). This
shows that $s=s'$, whence the contention.

(iii): In view of (ii), it suffices to show that the natural
map $F^+(X)\to\Hom_{\cC^\wedge}(h_\cS,F^+)$ is surjective for
every $\cS\in J(X)$. Hence, say that
$s\in\Hom_{\cC^\wedge}(h_\cS,F^+)$, and let $(S_i~|~i\in I)$
be a generating family for $\cS$.
By \eqref{subsec_interpret-descent}, $s$ corresponds to a
system $(s_i\in F^+(S_i)~|~i\in I)$ such that the following
holds. For every $i,j\in I$, and every pair of morphisms
$u_i:Y\to S_i$ and $u_j:Y\to S_j$ in $\cC\!/\!X$, we have
\set\begin{equation}\label{eq_repeat-plus}
F^+(u_i)(s_i)=F^+(u_j)(s_j).
\end{equation}
For every $i\in I$, let $\cS_i\in J(S_i)$ such that $s_i$
is the image of some $\bar s_i\in\Hom_{\cC^\wedge}(h_{\cS_i},F)$.
For every $u_i,u_j$ as above, set
$\cS_{ij}:=(\cS_i\times_{S_i}u_i)\cap(\cS_j\times_{S_j}u_j)$;
since $F$ is separated, \eqref{eq_repeat-plus} and claim
\ref{cl_better-injectivity} imply that the images of $\bar s_i$
and $\bar s_j$ agree in $\Hom_{\cC^\wedge}(h_{\cS_{ij}},F)$,
for every $i,j\in I$.

However, for every $i\in I$, let
$(g_{i\lambda}:T_{i\lambda}\to S_i~|~\lambda\in\Lambda_i)$ be a
generating family for $\cS_i$; then $\bar s_i$ corresponds to a
compatible system of sections $\bar s_{i\lambda}\in F(T_{i\lambda})$,
and $\cS_{ij}$ is the sieve of all morphisms $f:Z\to Y$ such that
$$
u_i\circ f=g_{i\lambda}\circ f'_i
\qquad\text{and}\qquad
u_j\circ f=g_{j\mu}\circ f'_j
$$
for some $\lambda\in\Lambda_i$, $\mu\in\Lambda_j$ and some
$f'_i:Z\to T_{i\lambda}$, $f'_j:Z\to T_{j\mu}$, so by construction
we have
\set\begin{equation}\label{eq_sorrento}
F(f'_i)(\bar s_{i\lambda})=F(f'_j)(\bar s_{j\mu})
\qquad
\text{for every $i,j\in I$ and $\lambda\in\Lambda_i$, $\mu\in\Lambda_\mu$}.
\end{equation}
Lastly, let $\cT$ be the sieve of $\cC\!/\!X$ generated by
$(g_i\circ g_{i\lambda}:T_{i\lambda}\to X~|~i\in I, \lambda\in\Lambda_i)$;
then $\cT$ covers $X$ (remark \ref{rem_topology}(iii)), and
\eqref{eq_sorrento} shows that the system
$(F(g_{i\lambda})(\bar s_{i\lambda})~|~i\in I, \lambda\in\Lambda_i)$
defines an element of $\Hom_{\cC^\wedge}(h_\cT,F)$, whose
image in $F^+(X)$ agrees with $s$.
\end{pfclaim}

From claim \ref{cl_plus-construction} we see that the rule :
$F\mapsto F^a:=(F^+)^+$ defines a left exact functor
$\cC^\wedge\to C^\sim$, with natural transformations
$\eta_F:F\Rightarrow i(F^a)$ for every $F\in\Ob(\cC^\wedge)$
and $\eps_G:(iG)^a\Rightarrow G$ for every $G\in\Ob(C^\sim)$
fulfilling the triangular identities of \eqref{subsec_adj-pair}.
The theorem follows.
\end{proof}

\begin{remark}\label{rem_rep-and-sheafify}
Let $C:=(\cC,J)$ be a small site.

(i)\ \
It has already been remarked that $C^\sim$ is complete, and
from theorem \ref{th_ass-sheaf} we also deduce that $C^\sim$
is cocomplete, and the functor \eqref{eq_asso-funct}
(resp. the inclusion functor $i:C^\sim\to\cC^\wedge$) commutes
with all colimits (resp. with all limits); more precisely,
if $F:I\to C^\sim$ is any functor from a small category $I$,
we have a natural isomorphism in $C^\sim$ (resp. in $\cC^\wedge$) :
$$
\colim_IF\isom(\colim_Ii\circ F)^a
\qquad
\text{(resp.\ \ \  $i(\lim_{I}F)\isom\lim_{I}i\circ F$)}
$$
(proposition \ref{prop_was-get-maddd}(iii,iv)); especially,
limits in $C^\sim$ are computed argumentwise. Also, all
colimits and all epimorphisms are universal in $C^\sim$ (see
examples \ref{ex_universal-col}(ii,iv) and
\ref{ex_strict-epis}(v,vii)), and filtered colimits in $C^\sim$
commute with finite limits, since the same holds in $\cC^\wedge$.

(ii)\ \ 
Furthermore, $C^\sim$ is well-powered and co-well-powered,
since the same holds for $\cC^\wedge$.
Especially, every morphism $f:F\to G$ in $C^\sim$ admits a
well defined image (see example \ref{ex_pre-misc}(i)) that
can be constructed explicitly as the subobject $(\Img\,i(f))^a$
(details left to the reader).

(iii)\ \
By composing with the Yoneda embedding, we obtain a functor
$$
h^a_C:\cC\to C^\sim \quad :\quad X\mapsto(h_X)^a
\qquad
\text{for every $X\in\Ob(\cC)$}
$$
and since the functor $(-)^a:\cC^\wedge\to C^\sim$ commutes
with small colimits, lemma \ref{lem_lable} yields a natural
isomorphism :
$$
\colim_{\cFib(F)}h^a_C\circ\ss_F\isom F
\qquad
\text{for every sheaf $F$ on $C$}.
$$

(iv)\ \
From the proof of claim \ref{cl_plus-construction} it is clear
that the functor
$$
\cC^\wedge\to C^\sep
\qquad
F\mapsto F^\sep:=\Img(F\to F^+)
$$
is left adjoint to the inclusion $C^\sep\to\cC^\wedge$.
Moreover, we have a natural identification :
$$
F^a\isom(F^\sep)^+
\qquad
\text{for every $F\in\Ob(\cC^\wedge)$}.
$$

(v)\ \ Let $\sV$ be a universe such that $\sU\subset\sV$;
from the definitions, it is clear that the fully faithful
inclusion $\cC^\wedge_\sU\subset\cC^\wedge_\sV$ restricts
to a fully faithful inclusion
$$
C^\sim_\sU\subset C^\sim_\sV.
$$
Moreover, by inspecting the proof of theorem \ref{th_ass-sheaf},
we deduce an essentially commutative diagram of categories :
$$
\xymatrix{ \cC^\wedge_\sU \ar[r] \ar[d] & C^\sim_\sU \ar[d] \\
           \cC^\wedge_\sV \ar[r] & C^\sim_\sV
}$$
whose vertical arrows are the inclusions, and whose horizontal
arrows are the functors $F\mapsto F^a$.
\end{remark}

In practice, one often encounters sites that are not small,
but which share many of the properties of small sites.
These more general situations are encompassed by the
following :

\begin{definition}\label{def_U-site}
Let $C:=(\cC,J)$ be a site.
\begin{enumerate}
\item
A {\em topologically generating family\/} for $C$ is a
subset $G\subset\Ob(\cC)$, such that, for every $X\in\Ob(\cC)$,
the family
$$
G/X:=\bigcup_{Y\in G}\Hom_\cC(Y,X)\subset\Ob(\cC\!/\!X)
$$
generates a sieve covering $X$.
\item
We say that $C$ is a {\em $\sU$-site}, if $\cC$ has small
$\Hom$-sets, and $C$ admits a small topologically generating
family. In this case, we also say that $J$ is a
{\em $\sU$-topology} on $\cC$.
\end{enumerate}
\end{definition}

\sset\subsubsection{}\label{subsec_U-site}
Let $C=:(\cC,J)$ be a $\sU$-site, and $G$ a small topologically
generating family for $C$. For every $X\in\Ob(\cC)$, denote by
$J_G(X)\subset J(X)$ the set of all sieves covering $X$ which
are generated by a subset of $G/X$ (notation of definition
\ref{def_U-site}(i)).

\begin{lemma}\label{lem_U-site}
With the notation of \eqref{subsec_U-site}, for every
$X\in\Ob(\cC)$ the following holds :
\begin{enumerate}
\item
$J_G(X)$ is a small set.
\item
$J_G(X)$ is a cofinal subset of $J(X)$ (for the partial
order given by inclusion of sieves).
\end{enumerate}
\end{lemma}
\begin{proof}(i) is left to the reader.

(ii): Let $\cS\in J(X)$, and say that $\cS$ is generated by
a family $(f_i:S_i\to X~|~i\in I)$ of objects of $\cC/X$
(indexed by some not necessarily small set $I$). Let $\cS'$
be the sieve generated by
$$
\bigcup_{i\in I}\{f_i\circ g~|~g\in G/S_i\}.
$$
It is easily seen that $\cS'\subset\cS$ and $\cS'\in J_G(X)$.
\end{proof}

\begin{remark}\label{rem_U-site}
(i)\ \ 
In the situation of \eqref{subsec_U-site}, let $\sV$ be a
universe with $\sU\subset\sV$, and such that $C$ is a
$\sV$-small site, so that we have a well defined functor
$(-)^+:\cC^\wedge_\sV\to\cC^\wedge_\sV$, as in \eqref{subsec_asso-topoi}.
Lemma \ref{lem_U-site} implies that the natural map :
$$
\colim_{\cS\in J_G(X)^o}\Hom_{\cC^\wedge}(h_\cS,F)\to F^+(X)
$$
is bijective. Therefore, if $F$ is a $\sU$-presheaf,
$F^+(X)$ is essentially $\sU$-small, and then the same
holds for $F^a(X)$. In other words, the restriction to
$\cC^\wedge_\sU$ of the functor $\cC^\wedge_\sV\to C^\sim_\sV$ :
$F\mapsto F^a$ is isomorphic to a functor that factors
through $C^\sim$. 

(ii)\ \
We deduce that theorem \ref{th_ass-sheaf} holds, more generally,
when $C$ is an arbitrary $\sU$-site. Likewise, a simple
inspection shows that remark \ref{rem_rep-and-sheafify}(i,iv,v)
holds when $C$ is only assumed to be a $\sU$-site.
\end{remark}

\begin{proposition}\label{prop_univ-effective}
Let $C:=(\cC,J)$ be a $\sU$-site. The following holds :
\begin{enumerate}
\item
A morphism in $C^\sim$ is an isomorphism if and only if it
is both a monomorphism and an epimorphism.
\item
All epimorphisms in $C^\sim$ are universal effective.
\end{enumerate}
\end{proposition}
\begin{proof}(i): Let $\phi:F\to G$ be a monomorphism in $C^\sim$;
then the morphism of presheaves $i(\phi):iF\to iG$ is also a
monomorphism, and it is easily seen that the cocartesian
diagram
$$
\cD \quad:\qquad
{\diagram iF \ar[r]^-{i(\phi)} \ar[d]_{i(\phi)} & iG \ar[d]^\alpha \\
          iG \ar[r] & iG\amalg_{iF}iG
\enddiagram}\qquad\qquad
$$
is also cartesian, hence the same holds for the induced diagram
of sheaves $\cD^a$. If moreover, $\phi$ is an epimorphism,
then $\alpha^a$ is an isomorphism, hence the same holds for
$\phi=(i(\phi))^a$.

(ii): Let $f:F\to G$ be an epimorphism in $C^\sim$; in view
of remarks \ref{rem_rep-and-sheafify}(i) and \ref{rem_U-site}(ii),
it suffices to show that $f$ is effective.
However, set $G':=\Img(i(f))$, and let $p_i:F\times_GF\to F$
(for $i=1,2$) be the two projections.
Suppose $\phi:F\to X$ is a morphism in $C^\sim$ such that
$\phi\circ p_1=\phi\circ p_2$; since
$i(F\times_GF)=iF\times_{iG}iF=iF\times_{G'}iF$, the morphism
$i(\phi)$ factors through a (unique) morphism $\psi:G'\to X$.
On the other hand, it is easily seen that $(G')^a=G$, hence
$\phi$ factors through the morphism $\psi^a:G\to X$.
\end{proof}

\begin{remark}\label{rem_mono-epi-fact}
Proposition \ref{prop_univ-effective}(i) implies that every
morphism $\phi:X\to Y$ in $C^\sim$ factors uniquely (up to
unique isomorphism) as the composition of an epimorphism
followed by a monomorphism. Indeed; such a factorization
is provided by the natural morphisms $X\to\Img(\phi)$ and
$\Img(\phi)\to Y$ (see example \ref{ex_pre-misc}(i)).
If $X\xrightarrow{\ \phi'\ } Z\to Y$ is another such
factorization, then by definition $\phi'$ factors through
a unique monomorphism $\psi:\Img(\phi)\to Z$.
However, $\psi$ is an epimorphism, since the same holds
for $\phi'$. Hence $\psi$ is an isomorphism.
\end{remark}

\begin{proposition}\label{prop_cover-is-iso}
Let $(\cC,J)$ be a $\sU$-site, $X\in\Ob(\cC)$, and $R$ any
subobject of $h_X$. The following conditions are equivalent :
\begin{enumerate}
\alphaenu
\item
The inclusion map
$i:R\to h_X$ induces an isomorphism on associated sheaves
$$
i^a:R^a\isom h_X^a.
$$
\item
$R$ covers $X$.
\end{enumerate}
\end{proposition}
\begin{proof}(b)$\Rightarrow$(a) : By definition, the natural
map $R^a(X)\to\Hom_{\cC^\wedge}(R,R^a)$ is bijective, hence there
exists a morphism $f:h_X\to R^a$ in $\cC^\wedge$ whose composition
with $i$ is the unit of adjunction $R\to R^a$.
Therefore, $f^a:h_X^a\to R^a$ is a left inverse for $i^a$.
On the other hand, we have a commutative diagram :
$$
\xymatrix{
\Hom_{\cC^\wedge}(h_X^a,h_X^a) \ar[r] \ar[d] &
\Hom_{\cC^\wedge}(R^a,h^a_X) \ar[d] \\
h_X^a(X) \ar[r] & \Hom_{\cC^\wedge}(R,h^a_X)
}$$
whose bottom and vertical arrows are bijective, so that the
same holds also for the top arrow. Set $g:=i^a\circ f^a$,
and notice that $g\circ i^a=i^a$, therefore $g$ must be
the identity of $h^a_X$, whence the contention.

(a)$\Rightarrow$(b): Let $\eta_X:h_X\to h_X^a$ be the unit of
adjunction, and set $j:=(i^a)^{-1}\circ\eta_X:h_X\to R^a$.
By remarks \ref{rem_rep-and-sheafify}(iv)  and \ref{rem_U-site}(ii),
we may find a covering subobject $i_1:R_1\to h_X$, and a morphism
$j_1:R_1\to R^\sep$ whose image in $\Hom_{\cC^\wedge}(R_1,R^a)$
equals $j\circ i_1$. Denote by $\eta'_X:h_X\to h^\sep_X$ the
unit of adjunction; by construction, the two morphisms
$$
i^\sep\circ j_1,\eta'_X\circ i_1:R_1\to h_X^\sep
$$
have the same image in $\Hom_{\cC^\wedge}(R_1,h^a_X)$.
This means that there exists a subobject $i_2:R_2\to R_1$
covering $X$, such that
$i^\sep\circ j_1\circ i_2=\eta'_X\circ i_1\circ i_2$.

Next, let $Y_\bullet:=(Y_\lambda\to X~|~\lambda\in\Lambda)$ be
a generating family for the sieve of $\cC/X$ corresponding to
$R_2$. There follows, for every $\lambda\in\Lambda$, a
commutative diagram :
\set\begin{equation}\label{eq_for-R-sep}
{\diagram
h_{Y_\lambda} \ar[r]^-{j_\lambda} \ar[d]_{i_\lambda} &
R^\sep \ar[d]^{i^\sep} \\
h_X \ar[r]^-{\eta_X'} & h_X^\sep.
\enddiagram}
\end{equation}
Then, for every $\lambda\in\Lambda$ there exists a covering
subobject $s_\lambda:R_\lambda\to h_{Y_\lambda}$ such that
$j_\lambda$ lifts to some $t_\lambda:R_\lambda\to R$, and we pick a
generating family $(Z_{\lambda\mu}\to Y_\lambda~|~\mu\in\Lambda_\lambda)$
for the sieve of $\cC/Y_\lambda$ corresponding to $R_\lambda$;
after replacing $Y_\bullet$ by the resulting family
$(Z_{\lambda\mu}\to X~|~\lambda\in\Lambda,\ \mu\in\Lambda_\lambda)$
(which still covers $X$, by virtue of remark \ref{rem_topology}(iii)),
we may assume that \eqref{eq_for-R-sep} lifts to a commutative diagram
$$
\xymatrix{
h_{Y_\lambda} \ar[r]^-{t_\lambda} \ar[d]_{i_\lambda} & R \ar[d] \\
h_X \ar[r]^-{\eta_X'} & h_X^\sep.}
$$
for every $\lambda\in\Lambda$. Then there exists a covering
subobject $s'_\lambda:R'_\lambda\to h_{Y_\lambda}$ such that
$i\circ t_\lambda\circ s'_\lambda=i_\lambda\circ s'_\lambda$
in $\Hom_{\cC^\wedge}(R'_\lambda,h_X)$. Lastly, set
$$
R':=\bigcup_{\lambda\in\Lambda}
\Img(i_\lambda\circ s'_\lambda:R'_\lambda\to h_X)
$$
(notice that $R'\in\Ob(\cC^\wedge_\sU)$ even in case $\Lambda$
is not a small set). It is easily seen that $R'$ is a covering
subobject of $X$, and the inclusion map $R'\to h_X$ factors
through $R$, so $R$ covers $X$ as well.
\end{proof}

\begin{definition}\label{def_tra-la-la}
Let $(\cC,J)$ be a site such that $\cC$ has small $\Hom$-sets,
and let $\phi:F\to G$ be a morphism in $\cC^\wedge$.
\begin{enumerate}
\item
We say that $\phi$ is a {\em covering morphism\/} if, for
every $X\in\Ob(\cC)$ and every morphism $h_X\to G$ in
$\cC^\wedge$, the image of the induced morphism
$F\times_Gh_X\to h_X$ is a covering subobject of $h_X$.
\item
We say that $\phi$ a {\em bicovering morphism\/} if both
$\phi$ and the diagonal morphism $F\to F\times_GF$ induced
by $\phi$, are covering morphisms.
\end{enumerate} 
\end{definition}

\begin{remark}\label{rem_iprippi}
(i)\ \
In the situation of definition \ref{def_tra-la-la}, let
$S:=\{X_i\to X~|~i\in I\}$ be a family of morphisms in $\cC$,
and pick a universe $\sV$ containing $\sU$, such that $I$
is $\sV$-small. Using remark \ref{rem_sieves-and-sub}(iii),
it is easily seen that $S$ covers $X$ if and only if the
induced morphism in $\cC^\wedge_\sV$
$$
\coprod_{i\in I}h_{X_i}\to h_X
$$
is a covering morphism.

(ii)\ \
With the notation of definition \ref{def_tra-la-la}, recall
that, by Yoneda's lemma, the set of morphisms $h_X\to G$ in
$\cC^\wedge$ is naturally identified with $GX$. Fix $s\in GX$,
and let $\psi_s:h_X\to G$ be the corresponding morphism; under
this identification, for every $Y\in\Ob(\cC)$ the image of
the induced map $(F\times_{(\phi,\psi_s)}h_X)(Y)\to h_X(Y)$ is
then the set of all morphisms $f:Y\to X$ such that $(Gf)(s)$
lies in the image of the map $\phi(f):FY\to GY$. Thus $\phi$
is a covering morphism if and only if for every $X\in\Ob(\cC)$
and every $s\in GX$, the sieve of $\cC/X$ generated by all
morphisms $f:Y\to X$ in $\cC$ such that
$(Gf)(s)\in\Img(\phi_Y:FY\to GY)$ covers $X$. In other words,
$\phi$ is a covering morphism if and only if for every
$X\in\Ob(\cC)$ and every $s\in GX$ there exists a covering family
$(f_i:Y_i\to X~|~i\in I)$ such that $(Gf_i)(s)\in\Img(\phi_{Y_i})$
for every $i\in I$.

(iii)\ \
Likewise, the diagonal morphism $F\to F\times_GF$ induced by
$\phi$ is a covering morphism if and only if for every
$X\in\Ob(\cC)$ and every two sections $s,s'\in FX$ such that
$\phi_X(s)=\phi_X(s')$, the set of all morphisms $f:Y\to X$
in $\cC$ such that $(Ff)(s)=(Ff)(s')$ generates a sieve covering
$X$. This is the same as saying that for every such $X$ and
every such pair $(s,s')$ there exists a covering family
$(f_i:Y_i\to X~|~i\in I)$ such that $(Ff_i)(s)=(Ff_i)(s')$
for every $i\in I$.
\end{remark}

\begin{corollary}\label{cor_bicover}
Let $C:=(\cC,J)$ be a $\sU$-site, and $\phi:F\to G$ a morphism
in $\cC^\wedge$. We have :
\begin{enumerate}
\item
$\phi$ is a covering morphism if and only if $\phi^a:F^a\to G^a$
is an epimorphism in $C^\sim$.
\item
The diagonal morphism $F\to F\times_GF$ induced by $\phi$ is
a covering morphism if and only if $\phi^a$ is a monomorphism
in $C^\sim$.
\item
$\phi$ is a bicovering morphism if and only if
$\phi^a:F^a\to G^a$ is an isomorphism in $C^\sim$.
\end{enumerate}
\end{corollary}
\begin{proof} Let $\sV$ be a universe with $\sU\subset\sV$,
and such that $\cC$ is $\sV$-small. Clearly $\phi$ is a
covering (resp. bicovering) morphism in $\cC_\sU^\wedge$ if
and only if the same holds for the image of $\phi$ under
the fully faithful inclusion $\cC^\wedge_\sU\subset\cC^\wedge_\sV$.
So we may replace $\sU$ by $\sV$, and assume that $C$ is a
small site. 

(i): Suppose that $\phi^a$ is an epimorphism; let $X$ be
any object of $\cC$, and $h_X\to G$ a morphism. Since the
epimorphisms of $C^\sim$ are universal (remark
\ref{rem_rep-and-sheafify}(i)), the induced morphism
$$
(\phi\times_GX)^a:(F\times_Gh_X)^a\to h^a_X
$$
is an epimorphism. Let $R\subset h_X$ be the image of
$\phi\times_Gh_X$; then the induced morphism $R^a\to h_X^a$
is both a monomorphism and an epimorphism, so it is an
isomorphism, by proposition \ref{prop_univ-effective}(i).
Hence $R$ is a covering subobject, according to proposition
\ref{prop_cover-is-iso}.

Conversely, suppose that $\phi$ is a covering morphism.
By remark \ref{rem_rep-and-sheafify}(iii), $G$ is the colimit
of a family $(h^a_{X_i}~|~i\in I)$ for certain $X_i\in\Ob(\cC)$.
By definition, the image $R_i$ of the induced morphism
$\phi\times_GX_i:F\times_Gh_{X_i}\to h_{X_i}$ covers $X_i$,
for every $i\in I$. Now, the induced morphism
$F\times_Gh_{X_i}\to R_i$ is an epimorphism, and
the morphism $R^a_i\to h^a_{X_i}$ is an isomorphism
(proposition \ref{prop_cover-is-iso}), hence $(\phi\times_GX_i)^a$
is an epimorphism, and then the same holds for
$$
\colim_{i\in I}(\phi\times_GX_i)^a:
\colim_{i\in I}F^a\times_{G^a}h^a_{X_i}\to\colim_{i\in I}h^a_{X_i}=G^a
$$
which is isomorphic to $\phi^a$, since the colimits of $C^\sim$
are universal (remark \ref{rem_rep-and-sheafify}(i)); so $\phi^a$
is an epimorphism.

(ii): If $\delta:F\to F\times_GF$ is a covering morphism, the
foregoing shows that the induced morphism
$\delta^a:F^a\to F^a\times_{G^a}F^a$ is both an epimorphism and
a monomorphism, hence it is an isomorphism (proposition
\ref{prop_univ-effective}(i)), so $\phi^a$ is a
monomorphism (proposition \ref{prop_diagonal-monos}(i)).
Conversely, if $\phi^a$ is a monomorphism, $\delta^a$ is an
isomorphism, hence $\delta$ is a covering morphism, by (i).

(iii) follows from (i), (ii) and proposition
\ref{prop_univ-effective}(i).
\end{proof}

\begin{remark}\label{rem_topol-on-presheaves}
(i)\ \
Let $C:=(\cC,J)$ be a site such that $\cC$ has small $\Hom$-sets,
and $\sV$ a universe containing $\sU$, such that $\cC^\wedge$ is
$\sV$-small. For every presheaf $F$ on $\cC$, let $J^\wedge(F)$ be
the set of all sieves $\cS\subset\cC^\wedge/F$ such that the natural
morphism
$$
\coprod_{(f:G\to F)\in\Ob(\cS)}\!\!\!\!\!\!G\to F
$$
is a covering morphism in $\cC^\wedge_\sV$. Notice that
$\cS\in J^\wedge(F)$ if and only if there exists a generating
family $(f_i:F_i\to F~|~i\in I)$ for $\cS$, indexed by a
$\sV$-small set $I$, such that the natural morphism
$\amalg_{i\in I}F_i\to F$ is a covering morphism. Since the
functor $(-)^a:\cC^\wedge_\sV\to C^\sim_\sV$ commutes with
arbitrary $\sV$-small colimits, corollary \ref{cor_bicover}(i)
says that the latter condition holds if and only if the induced
morphism $\amalg_{i\in I}F^a_i\to F^a$ is an epimorphism in
$C^\sim_\sV$.

(ii)\ \
We claim that the system $(J^\wedge(F)~|~F\in\Ob(\cC^\wedge))$
is a topology on $\cC^\wedge$. Indeed, it is clear from the
definition that if $\cS\in J^\wedge(F)$, then $\cS\times_Ff$
lies in $J^\wedge(F')$ for every morphism $f:F'\to F$ in
$\cC^\wedge$. Next, let $\cS,\cS'\subset\cC^\wedge/F$ be two
sieves, with $\cS\in J^\wedge(F)$ and such that
$\cS'\times_Ff\in J^\wedge(F')$ for every $(f:F'\to F)\in\Ob(\cS)$.
By the foregoing, this means that the induced morphism
$\amalg_{(g:G\to F)\in\Ob(\cS')}G^a\times_{F^a}F'^a\to F'^a$ is
an epimorphism in $C^\sim_\sV$ for every such $f$. Since
(in every category) an arbitrary coproduct of epimorphisms
is an epimorphism, it follows that both of the following
morphisms are epimorphisms :
$$
\coprod_{(f:F'\to F)\in\Ob(\cS)}\coprod_{(g:G\to F)\in\Ob(\cS')}
\!\!\!\!\!\!G^a\times_{F^a}F'^a\to\!\!\!\!\!\!\!\!\!\!\!\!
\coprod_{(f:F'\to F)\in\Ob(\cS)}\!\!\!\!\!\!F'^a\to F^a
$$
and then so is their composition. But 
$f\circ(g\times_FF'):G\times_FF'\to F$ lies in $\cS'$
for every $f\in\Ob(\cS)$ and every $g\in\Ob(\cS')$, so
$\cS'\in J^\wedge(F)$, as required. We denote the resulting
site
$$
C^\wedge:=(\cC^\wedge,J^\wedge).
$$
In light of remark \ref{rem_iprippi}(ii) it is easily seen
that $C^\wedge$ is independent of the choice of $\sV$.
\end{remark}

\subsection{Continuous and cocontinuous functors}
\label{sec_cont-functors}
We now begin the study of those functors that are compatible
with given topologies on their domain and codomain of definition;
as explained hereafter, such compatibility can manifest
itself in two distinct fashions :

\begin{definition}\label{def_dir-img-site}
Let $C=(\cC,J)$ and $C'=(\cC',J')$ be two sites, and
$g:\cC\to\cC'$ a functor. 

(i)\ \
We say that $g$ is {\em continuous\/} for the topologies $J$
and $J'$, if the following holds. For every universe $\sV$
and every $\sV$-sheaf $F$ on $C'$, the $\sV$-presheaf $g_\sV^\wedge F$
is a $\sV$-sheaf on $C$ (notation of \eqref{eq_pullback-presheaves}).
In this case, $g^\wedge_\sV$ clearly induces by restriction a functor
$$
\tilde g_{\sV*}:C^{\prime\sim}_\sV\to C^\sim_\sV.
$$

(ii)\ \
We say that $g$ is {\em cocontinuous\/} for the topologies
$J$ and $J'$ if the following holds. For every $X\in\Ob(\cC)$,
and every covering sieve $\cS'\in J'(gX)$, the sieve $g_{|X}^{-1}\cS'$
covers $X$ in $C$ (notation of definition \ref{def_sieve}(iii) and
\eqref{eq_restrict-over-X}).
\end{definition}

\begin{example}\label{ex_cont-map}
Let $T,T'$ be two topological spaces, and $f:T\to T'$
a continuous map; denote by $\cT$ (resp. $\cT'$) the
category of open subsets of $T$ (resp. of $\cT'$),
endowed with its natural topology as in example
\ref{ex_top-spaces}. Then $f$ induces a functor
$f^{-1}:\cT'\to\cT$ : $U\mapsto f^{-1}U$, and it is
easily seen that $f^*$ is a continuous functor,
according to definition \ref{def_dir-img-site}(i).
\end{example}

\begin{lemma}\label{lem_indepedence-V}
Let $C=(\cC,J)$ and $C'=(\cC',J')$ be two sites, $g:\cC\to\cC'$
a functor, and $\sV$ a universe such that $\cC$ is $\sV$-small
and $\cC'$ has small $\Hom$-sets. The following conditions are
equivalent :
\begin{enumerate}
\alphaenu
\item
For every $\sV$-sheaf $F$ on $C'$, the presheaf $g^\wedge_\sV F$
is a $\sV$-sheaf on $C$.
\item
For every morphism $\phi:F\to G$ of $\sV$-presheaves on $\cC$
that is bicovering for the topology $J$, the morphism
$g_{\sV!}(\phi):g_{\sV!}F\to g_{\sV!}G$ is bicovering for the
topology $J'$.
\item
For every $X\in\Ob(\cC)$ and every subobject $R\subset h_X$
covering $X$ for the topology $J$, the induced morphism of
presheaves $g_{\sV!}R\to g_{\sV!}h_X$ is bicovering for the
topology $J'$.
\item
$g$ is continuous for the topologies $J$ and $J'$.
\end{enumerate}
\end{lemma}
\begin{proof}(a)$\Rightarrow$(b): By virtue of remark
\ref{rem_was-cofinal}(iv) we may replace $\sV$ by a
larger universe, and assume that $\cC$ and $\cC'$ are
$\sV$-small. Then, by corollary \ref{cor_bicover}(iii),
the morphism $\phi^a:F^a\to G^a$ is an isomorphism, and
we need to check that $g_{\sV!}(\phi)^a$ is an isomorphism.
To this aim, let $H$ be any $\sV$-sheaf on $C'$; we get
a commutative diagram :
$$
\xymatrix@C-6pt{
\Hom_{\cC'^\wedge}(g_{\sV!}(G)^a,H) \ar[r] \ar[d] &
\Hom_{\cC'^\wedge}(g_{\sV!}(G),H) \ar[r] \ar[d] &
\Hom_{\cC^\wedge}(G,g^\wedge_\sV H) \ar[d] &
\Hom_{\cC^\wedge}(G^a,g^\wedge_\sV H) \ar[l] \ar[d] \\
\Hom_{\cC'^\wedge}(g_{\sV!}(F)^a,H) \ar[r] &
\Hom_{\cC'^\wedge}(g_{\sV!}(F),H) \ar[r] &
\Hom_{\cC^\wedge}(F,g^\wedge_\sV H) &
\Hom_{\cC^\wedge}(F^a,g^\wedge_\sV H) \ar[l]
}$$
whose horizontal arrows are bijections, due to (a).
Since $\phi^a$ is bijective, the right-most vertical arrow
is bijective, hence the same holds for the left-most vertical
arrow, whence the contention.

(b)$\Rightarrow$(c) and (d)$\Rightarrow$(a) are trivial.

(c)$\Rightarrow$(d): Let $\sV'$ be a universe, and $F$
a $\sV'$-sheaf on $C'$; we need to check that $g^\wedge_{\sV'}F$
is a sheaf on $C$. To this aim, we may replace $\sV'$ by
a larger universe, and assume that $\sV\subset\sV'$, and
that $\cC$ and $\cC'$ are $\sV'$-small; then by remark
\ref{rem_was-cofinal}(iv) and corollary \ref{cor_bicover}(iii)
we see that for every $X\in\Ob(\cC)$ and every covering subobject
$R\subset h_X$ the induced morphism $g_{\sV'!}(R)^a\to g_{\sV'!}(h_X)^a$
is an isomorphism. Now, we have a commutative diagram :
$$
\xymatrix{ \Hom_{\cC^\wedge}(h_X,g^\wedge_{\sV'}F) \ar[r] \ar[d] &
\Hom_{\cC'^\wedge}(g_{\sV'!}h_X,F) \ar[d] &
\Hom_{\cC'^\wedge}(g_{\sV'!}(h_X)^a,F) \ar[l] \ar[d] \\
\Hom_{\cC^\wedge}(R,g^\wedge_{\sV'}F) \ar[r] &
\Hom_{\cC'^\wedge}(g_{\sV'!}R,F) &
\Hom_{\cC'^\wedge}(g_{\sV'!}(R)^a,F) \ar[l]
}$$
whose horizontal arrows are bijections; moreover, the right-most
vertical arrow is bijective as well, hence the same holds for
the left-most vertical arrow, which is the contention.
\end{proof}

\begin{lemma}\label{lem_crit-continuity}
In the situation of definition {\em\ref{def_dir-img-site}},
consider the following conditions :
\begin{enumerate}
\alphaenu
\item
$g$ is continuous.
\item
For every covering family $(X_i\to X~|~i\in I)$ in $C$, the family
$(gX_i\to gX~|~i\in I)$ covers $gX$ in $C'$.
\item
For every small covering family $(X_i\to X~|~i\in I)$ in $C$,
the family $(gX_i\to gX~|~i\in I)$ covers $gX$ in $C'$.
\item
For every universe $\sV$, and every $\sV$-presheaf $F$ on
$\cC'$ that is separated for the topology $J'$, the presheaf
$g^\wedge_\sV F$ is separated for the topology $J$, so that
$g^\wedge_\sV$ restricts to a functor
$$
g^\sep_{\sV*}:C^{\prime\sep}_\sV\to C^\sep_\sV.
$$
\end{enumerate}
Then {\em (a)$\Rightarrow$(b)$\Leftrightarrow$(d)$\Rightarrow$(c)}.
Moreover, if all fibre products are representable in $\cC$, and
$g$ commutes with fibre products, then {\em (b)$\Rightarrow$(a)}.
Furthermore, if $C$ is a $\sU$-site, then {\em (c)$\Rightarrow$(b)}.
\end{lemma}
\begin{proof} Obviously (b)$\Rightarrow$(c).

(a)$\Rightarrow$(b): After replacing $\sU$ by a larger universe,
we may assume that $I$ is a small set and both $\cC$ and $\cC'$
are small. Let $F$ be any sheaf on $C'$; by assumption, $g^\wedge F$
is a sheaf on $C$, hence the natural map :
$$
F(gX)=g^\wedge F(X)\to
\prod_{i\in I}g^\wedge F(X_i)=\prod_{i\in I}F(gX_i)
$$
is injective (by \eqref{subsec_interpret-descent}).
This means that the induced morphism
\set\begin{equation}\label{eq_deduce-epi}
\coprod_{i\in I}h^a_{gX_i}\to h^a_{gX}
\end{equation}
is an epimorphism in $C^\sim$. Then the assertion follows from
corollary \ref{cor_bicover} and remark \ref{rem_iprippi}(i).

(d)$\Rightarrow$(b) is similar : we may assume that $I$, $\cC$
and $\cC'$ are small. If $F$ is any separated presheaf on $C'$,
then by assumption $g^\wedge F$ is separated on $C$, and arguing
as in the foregoing, we deduce that the induced morphism
$\amalg_{i\in I}h^\sep_{gX_i}\to h^\sep_{gX}$ is an epimorphism
in $C^\sep$, and then \eqref{eq_deduce-epi} is an epimorphism
in $C^\sim$, so we conclude, again by corollary \ref{cor_bicover}
and remark \ref{rem_iprippi}(i).

(b)$\Rightarrow$(d): let $F$ be a separated $\sV$-presheaf on
$C'$, and $(X_i\to X~|~i\in I)$ any covering family in $C$; in
view of \eqref{subsec_interpret-descent}, it suffices to show
that the induced map
$$
F(gX)=g^\wedge_\sV F(X)\to
\prod_{i\in I}g^\wedge_\sV F(X_i)=\prod_{i\in I}F(gX_i)
$$
is injective. But this is clear, since $(gX_i\to gX~|~i\in I)$
is a covering family in $C'$.

Next, suppose that (b) holds, the fibre products in
$\cC$ are representable, and $g$ commutes with all fibre
products. For every $i,j\in I$, set $X_{ij}:=X_i\times_XX_j$.
To show that (a) holds, it suffices -- in view of
\eqref{subsec_interpret-descent} -- to prove :

\begin{claim} The natural map
$$
g^\wedge F(X)\to\Hom_{\cC^\wedge}\Bigl(
\mathrm{Coequal}\Bigl( \coprod_{i,j\in I}
\xymatrix{\!h_{X_{ij}}\ar@<-.5ex>[r] \ar@<.5ex>[r] &} \!\!\!
\coprod_{i\in I}\!h_{X_i}\Bigr),g^\wedge F\Bigr)
$$
is bijective.
\end{claim}
\begin{pfclaim} Since $g_!$ is right exact, and due to
\eqref{eq_ess-comm-dig}, this is the same as the natural
map
$$
F(gX)\to\Hom_{\cC^\wedge}\Bigl(
\mathrm{Coequal}\Bigl( \coprod_{i,j\in I}
\xymatrix{\!h_{gX_{ij}}\ar@<-.5ex>[r] \ar@<.5ex>[r] &} \!\!\!
\coprod_{i\in I}\!h_{gX_i}\Bigr),F\Bigr).
$$
However, by assumption $gX_{ij}=gX_i\times_{gX}gX_j$,
and then the claim follows by applying
\eqref{subsec_interpret-descent} to the covering family
$(gX_i\to gX~|~i\in I)$.
\end{pfclaim}

Lastly, suppose that $C$ is a $\sU$-site; in order to show
that (c)$\Rightarrow$(b), we remark more precisely :

\begin{claim}\label{cl_darn-right}
Let $C$ be a $\sU$-site, $\cF:=(\phi_i:X_i\to X~|~i\in I)$ any
covering family. Then there exists a small set $J\subset I$
such that the subfamily $(\phi_i~|~i\in J)$ covers $X$.
\end{claim}
\begin{pfclaim} Let $\cS\subset\cC\!/\!X$ be the sieve
generated by $\cS$. By lemma \ref{lem_U-site}, we may find
a small covering family $\cF':=(\psi_i:X'_i\to X~|~i\in I')$
({\em i.e.} such that $I'$ is small), that generates a
sieve $\cS'\subset\cS$. Then, for every $i\in I'$ we
may find $\gamma(i)\in I$ such that $\psi_i$ factors
through $\phi_{\gamma(i)}$. The subset $J:=\gamma I'$
will do.
\end{pfclaim}

Let $\cF$ and $J$ be as in claim \ref{cl_darn-right}; then
the sieve $g\cS$ generated by $(g(\phi_i)~|~i\in I)$ contains
the sieve $g\cS'$ generated by $(g(\phi_i)~|~i\in J)$.
Especially, if $g\cS'$ is a covering sieve, the same holds
for $g\cS$, whence the contention.
\end{proof}

\sset\subsubsection{}\label{subsec_both-pigeons}
In the situation of definition \ref{def_dir-img-site}, let
$\sV$ be a universe with $\sU\subset\sV$, such that $C$ is
a $\sV$-site, and $\cC'$ has $\sV$-small $\Hom$-sets.
Then we may define a functor 
$$
\breve g{}_\sV^*:C_\sV^{\prime\sim}\to C_\sV^\sim
\qquad
F\mapsto(g^\wedge_\sV\circ i_{C'}F)^a
$$
(where $i_{C'}:C_\sV^{\prime\sim}\to\cC^\wedge_\sV$ is the forgetful
functor). On the other hand, if $C$ is a $\sV$-small site and
$C'$ is a $\sV$-site, we can define the functors
$$
g^a_{\sV!}:\cC^\wedge_\sV\to C'^\sim_\sV
\qquad
F\mapsto(g_{\sV!}F)^a
\qquad\text{and}\qquad
\tilde g^*_\sV:=g^a_{\sV!}\circ i_C:C^\sim_\sV\to C'^\sim_\sV.
$$
As usual, if $\sV=\sU$ we often omit the subscript $\sU$.
Later we shall generalize these constructions to the case
where both $C$ and $C'$ are $\sV$-sites : see corollary
\ref{cor_two-U-sites}. For now we remark :

\begin{lemma}\label{lem_breve}
In the situation of definition {\em\ref{def_dir-img-site}},
let $\sV$ be a universe with $\sU\subset\sV$, such that
$\cC$ is $\sV$-small, and $\cC'$ has $\sV$-small $\Hom$-sets.
Then the following conditions are equivalent :
\begin{enumerate}
\alphaenu
\item
$g$ is a cocontinuous functor.
\item
For every covering morphism $\phi:F\to G$ in $\cC'^\wedge$, the
morphism $g^\wedge(\phi):g^\wedge F\to g^\wedge G$ is covering.
\item
For every bicovering morphism $\phi:F\to G$ in $\cC'^\wedge$, the
morphism $g^\wedge(\phi):g^\wedge F\to g^\wedge G$ is bicovering.
\item
For every $X\in\Ob(\cC')$ and every covering subobject $R$ of
$h_X$, the induced morphism $g^\wedge R\to g^\wedge h_X$ is bicovering.
\item
For every $F\in\Ob(C^\sim_\sV)$, the $\sV$-presheaf $g_{\sV*}F$ is a
$\sV$-sheaf on $C'$ (see remark {\em\ref{rem_was-cofinal}(i)}).
\end{enumerate}
When these conditions hold, the restriction of $g_{\sV*}$ is a
right adjoint to $\breve g{}^*_\sV$ denoted
$$
\breve g_{\sV*}:C_\sV^\sim\to C_\sV^{\prime\sim}
$$
\end{lemma}
\begin{proof} After replacing $\sV$ by $\sU$, we may assume
that $\cC$ is small, and $\cC'$ has small $\Hom$-sets. Recall
that the unit of the natural adjunction for the pair $(g_!,g^\wedge)$
assigns to every $X\in\Ob(\cB)$ the morphism
$\eta_X:h_X\to g^\wedge(g_!h_X)\isom g^\wedge(h_{gX})$ such that
$\eta_X(s):=Fs$ for every $X'\in\Ob(\cC)$ and every
$(s:X'\to X)\in h_X(X')$ : see remark \ref{rem_was-cofinal}(v).
We notice :

\begin{claim}\label{cl_contends}
Let $X\in\Ob(\cC)$, and $\cT'\subset\cC'\!/\!gX$
a covering sieve; set $\cT:=g_{|X}^{-1}\cT'$. Then :
$$
h_\cT=g^\wedge(h_{\cT'})\times_{g^\wedge(h_{gX})}h_X.
$$
\end{claim}
\begin{pfclaim} For every $X'\in\Ob(\cC)$, the set
$(g^\wedge(h_{\cT'})\times_{g^\wedge(h_{gX})}h_X)(X')$ consists of
the pairs $(t,s)$ where $t:gX'\to gX$ is an object of $\cT'$
and $s:X'\to X$ is a morphism in $\cC$, such that $g(s)=t$.
In other words, this is the set of all $s\in h_X(X')$ such
that $g(s)\in\Ob(\cT')$. This is precisely the definition
of $h_\cT(X')$.
\end{pfclaim}

(a)$\Rightarrow$(b): Let $X\in\Ob(\cC)$ and
$\psi:h_X\to g^\wedge G$ any morphism in $\cC^\wedge$. By
adjunction, $\psi$ corresponds to a morphism $\psi':h_{gX}\to G$
in $\cC'^\wedge$, such that $\psi=g^\wedge(\psi')\circ\eta_X$.
By assumption, $\psi'$ induces a covering morphism
$F\times_Gh_{gX}\to h_{gX}$, whose image is therefore of the
form $h_{\cS'}$ for some covering sieve $\cS'\in J'(gX)$. Since
$g^\wedge$ is exact, the epimorphism $\pi:F\times_Gh_{gX}\to h_{\cS'}$
induces an epimorphism
$g^\wedge(\pi):g^\wedge(F\times_Gh_{gX}\to h_{gX})\to g^\wedge(h_{\cS'})$,
and notice the following commutative diagram whose three square
subdiagrams are cartesian :
$$
\xymatrix{ g^\wedge(F)\times_{g^\wedge(G)}h_X \ar[r] \ar[d] &
g^\wedge(h_{\cS'})\times_{g^\wedge(h_{gX})}h_X \ar[r] \ar[d] &
h_X \ar[d]^{\eta_X} \\
g^\wedge(F\times_Gh_{gX}) \ar[r]^-{g^\wedge(\pi)} \ar[d] &
g^\wedge(h_{\cS'}) \ar[r] & g^\wedge(h_{gX}) \ar[d]^{g^\wedge(\psi')} \\
g^\wedge(F) \ar[rr] & & g^\wedge(G)
}$$
In view of claim \ref{cl_contends}, we deduce that the image
of the morphism $g^\wedge(F)\times_{g^\wedge(G)}h_X\to h_X$ induced
by $g^\wedge(\phi)$ is $h_\cS$, where $\cS:=g^{-1}_{|X}\cS'$;
but $\cS$ is a covering sieve of $X$, since $g$ is cocontinuous.
This shows that (b) holds.

(b)$\Rightarrow$(c): If $\phi$ is a bicovering morphism, (b)
implies that both $g^\wedge(\phi)$ and the morphism
$\delta:g^\wedge(F)\to g^\wedge(F\times_GF)$ induced by $\phi$
are covering morphisms. However, $g^\wedge$ is exact, so $\delta$
is naturally identified with the diagonal morphism
$g^\wedge(F)\to g^\wedge(F)\times_{g^\wedge(G)}g^\wedge(F)$ induced
by $g^\wedge(\phi)$; thus, $g^\wedge(\phi)$ is bicovering.

(c)$\Rightarrow$(d) is obvious.

(d)$\Rightarrow$(e): Let $F$ be any sheaf on $C$; we have
to show that the natural map
$$
g_*F(Y)\to\Hom_{\cC^{\prime\wedge}}(R,g_*F)
$$
is bijective, for every $Y\in\Ob(\cC')$, and every covering 
subobject of $h_Y$. By adjunction (and by corollary
\ref{cor_bicover}), this is the same as saying that the
monomorphism $g^\wedge(R)\to g^\wedge(h_Y)$ is a covering
morphism, which holds by assumption (d).

(e)$\Rightarrow$(a): The assumption implies that, for every
sheaf $F$ on $C$, every $X\in\Ob(\cC)$, and every covering
sieve $\cS'\subset\cC'\!/\!gX$, the natural map
$$
g_*F(gX)\to\Hom_{\cC^{\prime\wedge}}(h_{\cS'},g_*F)
$$
is bijective. By adjunction, the same then holds for the natural
map
$$
\Hom_{\cC^\wedge}(g^\wedge(h_{gX}),F)\to
\Hom_{\cC^\wedge}(g^\wedge(h_{\cS'}),F)
$$
so the induced morphism $g^\wedge(h_{\cS'})\to g^\wedge(h_{gX})$
is bicovering (proposition \ref{prop_cover-is-iso}). Also,
this morphism is a monomorphism (since $g^\wedge$ commutes
with all limits); therefore, after base change along the
unit of adjunction $\eta_X:h_X\to g^\wedge(h_{gX})$, we deduce
a covering monomorphism
$$
g^\wedge(h_{\cS'})\times_{g^\wedge(h_{gX})}h_X\to h_X.
$$
Then the contention follows from claim \ref{cl_contends}.

Lastly, the assertion concerning the left adjoint
$\breve g{}_\sU^*$ is immediate from the definitions.
\end{proof}

\begin{lemma}\label{lem_cont-funct-site}
Let $C':=(\cC',J')$ be a $\sU$-site, $C:=(\cC,J)$ a small site,
and $g:\cC\to\cC'$ a continuous functor. Then the following
holds :
\begin{enumerate}
\item
The functor $\tilde g{}_\sU^*:C_\sU^\sim\to C_\sU^{\prime\sim}$ of
\eqref{subsec_both-pigeons} is left adjoint to $\tilde g_{\sU*}$.
\item
The natural diagrams of functors :
$$
\xymatrix{
\cC \ar[r]^-g \ar[d]_{h^a_\cC} & \cC' \ar[d]^{h^a_{\cC'}} &
\cC^\wedge_\sU \ar[rr]^-{(-)^a} \ar[rd]_{g^a_{\sU!}} & &
C^\sim_\sU \ar[ld]^{\tilde g{}^*_\sU} \\
C^\sim_\sU \ar[r]^-{\tilde g{}_\sU^*} & C_\sU^{\prime\sim}
& & C'^\sim_\sU
}$$
are essentially commutative. (Notation of theorem
{\em\ref{th_ass-sheaf}}.)
\end{enumerate}
\end{lemma}
\begin{proof} The first assertion is straightforward,
and the commutativity of the square diagram in (ii) is
reduced to the corresponding assertion for
\eqref{eq_ess-comm-dig}, which has already been remarked.
To show the commutativity for the triangular diagram in (ii),
it suffices to check that for every presheaf $F$ on $\cC$
and every sheaf $G$ of $C'$ the natural map
$\eta_F:F\to i_C(F^a)$ induces a bijection
$$
\Hom_{C'^\sim}(\tilde g{}^*_\sU(F^a),G)\isom
\Hom_{C'^\sim}(g^a_{\sU!}(F),G).
$$
However, by adjunction, the latter is naturally identified
with the map
\set\begin{equation}\label{eq_induced-by-eta_F}
\Hom_{\cC^\wedge}(F^a,g^\wedge_\sU G)\to\Hom_{\cC^\wedge}(F,g^\wedge_\sU G)
\end{equation}
induced by $\eta_F$; but $g^\wedge_\sU G$ is a sheaf, since $g$
is continuous, so \eqref{eq_induced-by-eta_F} is indeed
bijective.
\end{proof}

\sset\subsubsection{}\label{subsec_adj-universes}
Let $(\cC,J)$ be a small site, $(\cC',J')$ a $\sU$-site,
$u:\cC\to\cC'$ a continuous (resp. cocontinuous) functor.
Let also $\sV$ be a universe such that $\sU\subset\sV$.
Then it follows easily from remark \eqref{rem_was-cofinal}(iv)
and lemma \ref{lem_cont-funct-site}(i) (resp. from lemma
\ref{lem_breve}), and from remark \ref{rem_U-site}(ii) that
we have essentially commutative diagrams of categories :
$$
{\diagram C^\sim_\sU \ar@<.5ex>[r]^-{\tilde u{}_\sU^*} \ar[d] &
\ar@<.5ex>[l]^-{\tilde u{}_{\sU*}} C^{\prime\sim}_\sU \ar[d] \\
C^\sim_\sV \ar@<.5ex>[r]^-{\tilde u{}_\sV^*} &
\ar@<.5ex>[l]^-{\tilde u{}_{\sV*}} C^{\prime\sim}_\sV
\enddiagram}
\qquad\text{(resp.}\qquad
{\diagram C^\sim_\sU \ar@<.5ex>[r]^-{\breve u{}_{\sU*}} \ar[d] &
\ar@<.5ex>[l]^-{\breve u{}_\sU^*} C^{\prime\sim}_\sU \ar[d]\\
C^\sim_\sV \ar@<.5ex>[r]^-{\breve u{}_{\sV*}} & C^{\prime\sim}_\sV
\ar@<.5ex>[l]^-{\breve u{}_\sV^*}
\enddiagram}\text{)}$$
whose vertical arrows are the inclusion functors. More
generally, the diagram for $\breve u_{\sU*}$ is well defined
and essentially commutative, whenever $C$ is a $\sU$-site,
and $\cC'$ has small $\Hom$-sets.

\begin{lemma}\label{lem_needed}
{\em (i)}\ \
Let $C:=(\cC,J)$ and $C':=(\cC',J')$ be two sites, and
$v:\cC\to\cC'$, $u:\cC'\to\cC$ two functors, such that
$v$ is left adjoint to $u$. The following conditions
are equivalent:
\begin{enumerate}
\alphaenu
\item
$u$ is continuous.
\item
$v$ is cocontinuous.
\end{enumerate}

{\em (ii)}\ \
Moreover, when these conditions hold, then for every universe
$\sV$ such that $\cC$ and $\cC'$ are $\sV$-small, we have
natural isomorphisms of functors :
$$
\tilde u_{\sV*}\isom\breve v_{\sV*}
\qquad
\tilde u{}_\sV^*\isom\breve v{}_\sV^*.
$$
\end{lemma}
\begin{proof} In view of lemma \ref{lem_indepedence-V}, we
may replace $\sU$ by a larger universe, after which we may
assume that $C$ and $C'$ are small sites. In this case, the
lemma follows from lemma \ref{lem_breve} and proposition
\ref{prop_was-get-maddd}(i).
\end{proof}

\begin{lemma}\label{lem_improve}
Let $(\cC,J)$ be a small site, $(\cC',J')$ a $\sU$-site,
$u:\cC\to\cC'$ a continuous and cocontinuous functor.
Then we have :
\begin{enumerate}
\item
$\tilde u_*=\breve u{}^*$ and this functor admits the left
adjoint $\tilde u{}^*$ and the right adjoint $\breve u_*$.
\item
We have a natural isomorphism of functors
$$
(-)^a\circ u^\wedge\isom\tilde u_*\circ(-)^a
\ :\ \cC'^\wedge\to C^\sim.
$$
\item
$\tilde u{}^*$ is fully faithful if and only if the same holds
for $\breve u_*$.
\item
If $u$ is fully faithful, then the same holds for $\tilde u{}^*$.
The converse holds, provided the topologies $J$ and $J'$ are
coarser than the canonical topologies. 
\end{enumerate}
\end{lemma}
\begin{proof}(i) is clear by inspecting the definitions.

(ii): For every presheaf $G$ on $\cC'$, let $i_G:G\to G^a$
be the natural morphism; due to lemma \ref{lem_breve}
the induced morphism $u^\wedge(i_G):u^\wedge G\to u^\wedge(G^a)$
is bicovering, and notice that $u^\wedge(G^a)$ is a sheaf,
since $u$ is continuous. Thus, $u^\wedge(i_G)$ factors through
the natural morphism $i_{u^\wedge G}:u^\wedge G\to(u^\wedge G)^a$
and a unique morphism in $C^\sim$
$$
\omega_G:(u^\wedge G)^a\to u^\wedge(G^a).
$$
Since $i_{u^\wedge G}$ is bicovering, the same holds for
$\omega_G$, {\em i.e.} the latter is an isomorphism of
on $C^\sim$. It is then easily seen that the rule
$G\mapsto\omega_G$ yields the sought isomorphism.

Assertion (iii) follows from (i) and proposition
\ref{prop_fullfaith-adjts}(iv). Next, suppose that $u$
is fully faithful; then the same holds for $\breve u_*$
(corollary \ref{cor_lable}(ii)), so the claim follows from
(iii). Lastly, suppose that $\tilde u{}^*$ is fully faithful,
and both $J$ and $J'$ are coarser than the
canonical topologies on $\cC$ and $\cC'$. In such case,
the Yoneda embedding for $\cC$ (resp. $\cC'$) realizes
$\cC$ (resp. $\cC'$) as a full subcategory of $C^\sim$
(resp. of $C^{\prime\sim}$). Then, from \eqref{eq_ess-comm-dig}
and the explicit expression of $\tilde u{}^*$ provided by
lemma \ref{lem_cont-funct-site}(i), we deduce that $u$
is fully faithful.
\end{proof}

\begin{definition}\label{def_induced-top}
Let $C:=(\cC,J)$ be a site, $\cB$ any category, and $g:\cB\to\cC$
a functor. Pick a universe $\sV$ such that $\cB$ is $\sV$-small
and $C$ is a $\sV$-site. According to remark \ref{rem_sheaves}(iv),
there is a finest topology $J_g$ on $\cB$ such that, for every
$\sV$-sheaf $F$ on $C$, the $\sV$-presheaf $g^\wedge F$ is a $\sV$-sheaf
on $(\cB,J_g)$. By lemma \ref{lem_indepedence-V}, the topology
$J_g$ is independent of the chosen universe $\sV$. We call $J_g$
the {\em topology induced via $g$ by $J$ on $\cB$}. Clearly $g$ is
continuous for the sites $(\cB,J_g)$ and $C$.
\end{definition}

\begin{lemma}\label{lem_induced-top}
In the situation of definition {\em\ref{def_induced-top}}, let
$X$ be any object of $\cB$, and $R\subset h_X$ any
subobject in $\cB^\wedge$. We have :
\begin{enumerate}
\item
$R\in J_g(X)$ if and only if for every morphism $Y\to X$ in
$\cB$, the induced morphism $g_{\sV!}(R\times_XY)\to h_{g(Y)}$
is a bicovering morphism in $\cC^\wedge_\sV$ (for the topology
$J$ on $\cC$).
\item
Suppose moreover, that either one of the following condition holds:
\begin{enumerate}
\item
The functor $g^a_{\sV!}:\cB^\wedge_\sV\to C^\sim_\sV$ commutes with all
fibre products.
\item
All fibre products are representable in $\cB$, and $g$ commutes
with fibre products.
\end{enumerate}
Then a family $(f_i:B_i\to B~|~i\in I)$ of morphisms of $\cB$
generates a covering sieve of $J_g$ if and only if
$(g(f_i):gB_i\to gB~|~i\in I)$ generates a covering sieve of $J$.
\end{enumerate}
\end{lemma}
\begin{proof}(i): According to remark \ref{rem_sheaves}(iii),
we have $R\in J_g(X)$ if and only if for every morphism
$Y\to X$ in $\cC$ the natural morphism
$\phi_Y:R\times_XY\to h_Y$ induces a bijection
$$
(g^\wedge_\sV G)(Y)\isom\Hom_{\cB^\wedge}(R\times_XY,g^\wedge_\sV G)
\qquad
\text{for every $G\in\Ob(C^\sim_\sV)$}.
$$
By adjunction, this is equivalent to saying that $\phi_Y$
induces a bijection
$$
G(gY)\isom\Hom_{\cC^\wedge}(g_{\sV!}(R\times_XY),G)
\qquad
\text{for every $G\in\Ob(C^\sim_\sV)$}.
$$
The latter is in turn naturally identified with the map
$$
\Hom_{C^\sim_\sV}(g_{\sV!}(\phi_Y)^a,G):
\Hom_{C^\sim_\sV}(h^a_{gY},G)\to\Hom_{C^\sim_\sV}(g_{\sV!}(R\times_XY)^a,G).
$$
By Yoneda's lemma, it follows that $R\in J_g(R)$ if and only if
$g_{\sV!}(\phi_Y)^a$ is an isomorphism, and this is equivalent
to the stated condition, by corollary \ref{cor_bicover}.

(ii): Suppose first that condition (a) holds. Notice that
since $\cB$ is $\sV$-small, we may assume that the indexing
set $I$ is $\sV$-small as well. Let $R\subset h_B$ (resp.
$R'\subset h_{gB}$) be the subobject generated by the family
$(f_i~|~i\in I)$ (resp. $(g(f_i)~|~i\in I)$); if $R$ is a
covering subobject of $h_B$ for $J_g$, then $R'$ is a covering
subobject of $h_{gB}$ for $J$, by lemma \ref{lem_crit-continuity}.
Conversely, suppose that $R'$ is a covering subobject of $h_{gB}$
for $J$; according to (i) and corollary \ref{cor_bicover}, we
need to check that the inclusion $j:R\to h_B$ induces an
isomorphism $g^a_{\sV!}(R\times_BB')\to h^a_{gB'}$ for every
morphism $B'\to B$ in $\cB$. Since $g^a_{\sV!}$ commutes with
fibre products, we are then easily reduced to checking that
$j$ induces an isomorphism $g^a_{\sV!}(j):g^a_{\sV!}R\to h^a_B$.
However, $g^a_{\sV!}(j)$ is a monomorphism, since $g^a_{\sV!}$
commutes with fibre products (proposition
\ref{prop_was-also-cofinal}(i)). To show that $g^a_{\sV!}(j)$ is an
epimorphism, let $S:=\coprod_{i\in I}h_{B_i}$ (which is an object
of $\cB^\wedge_\sV$, since $I$ is $\sV$-small) and denote by
$j':S\to h_B$ the morphism induced by the morphisms $f_i$.
The image of $g^a_{\sV!}(j)$ contains the image of
$g^a_{\sV!}(j')$; on the other hand, $g^a_{\sV!}$ commutes with
coproducts, since it is a left adjoint, and we have natural
identifications : $g^a_{\sV!}(h_{B_i})=h^a_{gB_i}$ for every $i\in I$.
Thus, the image of $g^a_{\sV!}(j')$ equals that of the morphism
$\phi:\coprod_{i\in I}h^a_{gB_i}\to h^a_B$ induced by the morphisms
$g(f_i)$. But since $R'$ is a covering subobject, $\phi$ is an
epimorphism, so the same holds for $g^a_{\sV!}(j)$.

Next, suppose that condition (b) holds, and let
$f_\bullet:=(f_i:B_i\to B~|~i\in I)$ be a family of morphisms of
$\cB$ such that $(g(f_i)~|~i\in I)$ covers $gB$ in the topology
$J$. Set $B_{ij}:=B_i\times_BB_j$ for every $i,j\in I$; the
subobject $R\subset h_B$ generated by $f_\bullet$ is the
coequalizer in $\cB^\wedge$ of the natural projections
$$
\xymatrix{
\coprod_{i,j\in I} h_{B_{ij}} \ar@<.5ex>[r] \ar@<-.5ex>[r] &
\coprod_{i\in I} h_{B_i}
}$$
and according to (i) we need to check that for every morphism
$B'\to B$ the induced morphism $g_{\sV!}(R\times_BB')\to h_{gB'}$
is a bicovering morphism. However, set as well
$B'_i:=B_i\times_BB'$ and $B'_{ij}:=B_{ij}\times_BB'$ for every
$i,j\in I$; then $R\times_BB'$ is the coequalizer of the
induced projections
$$
\xymatrix{
\coprod_{i,j\in I} h_{B'_{ij}} \ar@<.5ex>[r] \ar@<-.5ex>[r] &
\coprod_{i\in I} h_{B'_i}.
}$$
On the other hand, $g_{\sV!}$ commutes with coequalizers and
coproducts, since it is a left adjoint; in view of remark
\ref{rem_was-cofinal}(ii) we deduce that
$g_{\sV!}(R\times_BB')$ is the coequalizer of the projections
$$
\xymatrix{
\coprod_{i,j\in I} h_{gB'_{ij}} \ar@<.5ex>[r] \ar@<-.5ex>[r] &
\coprod_{i\in I} h_{gB'_i}.
}$$
But by assumption we have as well $gB'_i=gB_i\times_{gB}gB'$,
whence $h_{gB_i}=h_{gB'_i}\times_{h_{gB}}h_{B'}$, and likewise
for $h_{gB'_{ij}}$, for every $i,j\in I$ (corollary
\ref{cor_pre-misc}(vi)). We conclude that $g_{\sV!}(R\times_BB')$
is the subobject of $h_{B'}$ generated by the family
$(g(f_i)\times_{gB}gB'\to gB'~|~i\in I)$, and it is therefore
a covering subobject, as required.
\end{proof}

\begin{remark}\label{rem_change-universe}
(i)\ \
In the situation of lemma \ref{lem_induced-top}, let $\sV'$
be another universe with $\sV\subset\sV'$, and $F$ a
$\sV'$-presheaf on $\cB$. According to lemma \ref{lem_lable},
$F$ is naturally isomorphic to the colimit of the functor
$h_\cB\circ\ss_F:\cFib(F)\to\cB^\wedge_\sV\subset\cB^\wedge_{\sV'}$,
where $\ss_F:\cFib(F)\to\cB$ is the source functor. For every
finite subset $S\subset\Ob(\cFib(F))$, let $F_S\subset F$ be
the subobject generated by the union of the images of all the
morphisms $h_{(X,s)}:h_B\to F$ with $(B,s)\in S$ (notation of
\eqref{subsec_category-of-elements}), and notice that $F_S$ is
isomorphic to an object of $\cB^\wedge_\sV$. After replacing each
$F_S$ by an isomorphic object of $\cB^\wedge_\sV$, we obtain therefore
$F$ as the filtered colimit of a system $(F_S~|~S\in\Sigma_F)$
of objects of $\cB^\wedge_\sV$, indexed by the set $\Sigma_F$
of all finite subsets of $\Ob(\cFib(F))$.

(ii)\ \
Now, let $F_1\xrightarrow{\phi_1}F_0\xleftarrow{\phi_2}F_2$ be
two morphisms in $\cB^\wedge_\sV$. Consider the set
$\Sigma(\phi_1,\phi_2)$ of all triples $(S_0,S_1,S_2)$ where
$S_i$ is a finite subset of $\Ob(\cFib(F_i))$ for $i=0,1,2$, and
$(B,\phi_j\circ s)\in S_0$ for $j=1,2$ and every $(B,s)\in S_j$.
We endow $\Sigma(\phi_1,\phi_2)$ with the partial order such
that $(S_0,S_1,S_2)\geq(S'_0,S'_1,S'_2)$ if and only if
$S'_i\subset S_i$ for $i=0,1,2$. We have a natural functor
$$
\tau:\Sigma(\phi_1,\phi_2)\to\cB^\wedge_\sV
\qquad
(S_0,S_1,S_2)\mapsto F_{S_1}\times_{F_{S_0}}F_{S_2}
$$
and by a simple inspection we see that the colimit of $\tau$
is naturally isomorphic to $F_1\times_{F_0}F_2$.

(iii)\ \
We deduce that $G:=g^a_{\sV'!}(F_1\times_{F_0}F_2)$ is naturally
isomorphic to the colimit of $g^a_{\sV'!}\circ\tau$, since
$g^a_{\sV'!}$ is a left adjoint; but $g^a_{\sV'!}\circ\tau$ is
in turn naturally isomorphic to the functor $g^a_{\sV!}\circ\tau$,
by lemma \ref{lem_cont-funct-site}(ii). Suppose now that
$g^a_{\sV!}$ commutes with all fibre products of $\cB^\wedge_\sV$.
Then $G$ is also naturally isomorphic to the filtered colimit
of the system
$$
(g^a_{\sV!}F_{1,S_1}\times_{g^a_{\sV!}F_{0,S_0}}g^a_{\sV!}F_{2,S_2}
~|~(S_0,S_1,S_2)\in\Sigma(\phi_1,\phi_2)).
$$
But since filtered colimits in $\cB^\wedge_{\sV'}$ commute with
fibre products, the latter is naturally isomorphic to the fibre
product $G_1\times_{G_0}G_2$, where $G_i$ is the filtered colimit
of the system $(g^a_{\sV'!}F_{S_i}~|~S_i\in\Sigma_{F_i})$ (notice
that the projections $\Sigma(\phi_1,\phi_2)\to\Sigma_{F_i}$ :
$(S_0,S_1,S_2)\to S_i$ are obviously final functors). But again,
since $g^a_{\sV'!}$ commutes with filtered colimits, $G_i$ is
naturally isomorphic to $g^a_{\sV'!}(F_i)$. We
conclude that $g^a_{\sV'!}$ commutes with fibre products as
well. In particular, condition (ii.a) of lemma
\ref{lem_induced-top} is independent of the choice of
universe $\sV$.
\end{remark}

\subsection{Morphisms of sites}
We observed in example \ref{ex_cont-map} that every continuous
map of topological spaces $f:T\to T'$ induces a continuous
functor $u:=f^{-1}:\cT'\to\cT$ between the corresponding
categories of open subsets.
It is well known that the induced functor on categories
of sheaves $\tilde u{}^*:\cT^\sim\to\cT'^\sim$ is moreover
{\em exact}, whereas for a general continuous functor $g$
between arbitrary sites, only the right exactness of
$\tilde g{}^*$ is always assured. The exactness of $\tilde g{}^*$
therefore singles out an interesting class of continuous
functors, which will play an important role in our discussion
of topoi, and -- even more crucially --  for the study of
the functorial properties of the categories of stacks,
in section \ref{sec_funct-stacks}. In this section we carry
out a preliminary investigation of this class of functors,
and prove a useful characterization. Let us begin with :

\begin{definition}\label{def_morph-of-sites}
(i)\ \
Let $C=(\cC,J)$ and $C'=(\cC',J')$ be two sites, and
$\sV$ a universe such that $C$ and $C'$ are $\sV$-small.
A {\em morphism of sites\/} $C'\to C$ is the datum of
a continuous functor $g:\cC\to\cC'$, such that the left
adjoint $\tilde g{}^*_\sV$ of $\tilde g_{\sV*}$ is exact
(notation of definition \ref{def_dir-img-site}(i)).

(ii)\ \
Let $\sU'$ be a universe with $\sU\in\sU'$; the $\sU'$-small
$\sU$-sites are the objects of a $2$-category
$$
(\sU,\sU')\tdu\Site
$$
whose $1$-cells are the morphisms of sites, and whose $2$-cells
$\beta:g\Rightarrow g'$ are the natural transformations
$g'\Rightarrow g$ between the underlying functors of such
morphisms, with the obvious composition laws for $1$-cells
and $2$-cells.

(iii)\ \
We shall also be interested in the $2$-category of
{\em $\sU$-lex-sites}
$$
(\sU,\sU')\tdu\lex.\Site
$$
defined as the $2$-subcategory of $(\sU,\sU')\tdu\Site$ whose
objects are the finitely complete $\sU'$-small $\sU$-sites and
whose $1$-cells are the morphisms of sites $g$ as in (i) such
that the underlying functor $g:\cC\to\cC'$ is left exact. For
any two such morphisms $g$, $g'$, the $2$-cells $g\Rightarrow g'$
are the same as in $(\sU,\sU')\tdu\Site$.
\end{definition}

\begin{remark}\label{rem_morph-of-sites}
(i)\ \
With the notation of definition \ref{def_morph-of-sites}(ii),
two different choices of the auxiliary universe $\sU'$
{\em do not} necessarily yield $2$-equivalent $2$-categories;
nevertheless, such choices are usually ininfluent in the proofs
of our results, so we mostly omit mentioning them explicitly,
and write $\sU\tdu\Site$ and $\sU\tdu\lex.\Site$ for these
$2$-categories. When the choice of $\sU$ is clear from the
context, we shall likewise drop the mention of $\sU$, and
write just $\Site$ and $\lex.\Site$.

(ii)\ \
On the other hand, the following proposition
\ref{prop_morph-of-sites} will show that the definition
of morphism of sites $g:(\cC',J')\to(\cC,J)$ depends only
on the underlying functor $g:\cC\to\cC'$, and not on the
choice of a universe $\sV$ such that $C$ and $C'$ are
$\sV$-small.

(iii)\ \
Let $C:=(\cC,J)$ and $C':=(\cC',J')$ be two sites, $u:C'\to C$
a morphism of sites, and $S$ a set. The {\em constant presheaf
on $\cC$ with value $S$} is the constant functor $c_S:\cC^o\to\Set$
associated with $S$ (see \eqref{subsec_wishful}), and the
{\em constant sheaf on $C$ with value $S$} is the associated
sheaf $c^a_S$. Clearly we have $c_S=\amalg_{s\in S}c_{\{s\}}$,
and each presheaf $c_{\{s\}}$ is a final object of the category
$\cC^\wedge$. It follows that $c^a_S=\amalg_{s\in S}c^a_{\{s\}}$,
and $c^a_{\{s\}}$ is a final object of $C^\sim$ (example
\ref{ex_equalizers}(iv)). Since $\tilde u{}^*$ commutes with
all colimits and is left exact by assumption, we deduce that
$\tilde u{}^*(c^a_S)$ is naturally isomorphic to the constant
sheaf with value $S$ on $C'$. This assertion may fail for
general continuous functors.
\end{remark}

For our characterization of morphisms of sites, we shall
need the following :

\begin{definition}\label{def_locally-cofiltered}
Let $(\cB,J)$ be a site, $F:\cC\to\cB$ a fibration,
and fix a cleavage $\blambda$ for $F$, so that for every
$X\in\Ob(\cC)$ and every morphism $f:B\to FX$ we have the
cartesian morphism $\blambda(X,f):f^*X\to X$. We say that
$F$ is {\em locally cofiltered} (relative to the topology $J$)
if the following holds for every $B\in\Ob(\cB)$ and every
$X,X'\in\Ob(\cC_B)$ :
\begin{enumerate}
\alphaenu
\item
There exists a covering family $(B_i\to B~|~i\in I)$ in the
topology $J$ with $\Ob(\cC_{B_i})\neq\emptyset$ for every
$i\in I$.
\item
There exist a covering family $(f_i:B_i\to B~|~i\in I)$ in the
topology $J$, and for every $i\in I$ an object $Y_i\in\Ob(\cC_{B_i})$
with two morphisms $f^*_iX\leftarrow Y_i\to f^*_iX'$ in $\cC_{B_i}$.
\item
For every pair of morphisms $g,g':X'\to X$ in $\cC_B$ there exist
a covering family $(f_i:B_i\to B~|~i\in I)$ in the topology $J$, and
for every $i\in I$ a morphism $h_i:Y_i\to f^*_iX'$ in $\cC_{B_i}$ such
that $f^*_i(g)\circ h_i=f^*_i(g')\circ h_i$.
\end{enumerate}
\end{definition}

\sset\subsubsection{}\label{subsec_Fubini-sheaf}
In the situation of definition \ref{def_locally-cofiltered}, suppose
that $(\cB,J)$ and $\cC$ are small; we consider the functor
$$
\int^\blambda_a:\cC^\wedge\to(\cB,J)^\sim
\qquad
G\mapsto\Bigl(\int^\blambda G\Bigr)^a
$$
where $\int^\blambda:\cC^\wedge\to\cB^\wedge$ is defined as in
example \ref{ex_Fubini-fibred}(iv). Recall that the sections
of $\int^\blambda G(B)$ are the equivalence classes $[X,t]$
of pairs $(X,t)$ consisting of an object $X$ of $F^{-1}B$
and a section $t\in GX$, for the equivalence relations
on the set of such pairs explicited in example
\ref{ex_complete-cats}(i).

\begin{lemma}\label{lem_locally-filt-equality}
With the notation of \eqref{subsec_Fubini-sheaf}, let $L$
be any presheaf on $\cC$, and $M:=\int^\blambda L$. Let also
$B\in\Ob(\cB)$ and $(X,t),(X',t')$ two pairs consisting of
objects $X,X'\in F^{-1}B$ and sections $t\in LX$, $t'\in LX'$.
If $[X,t]=[X',t']$ in $M(B)$, there exist a covering family
$(f_j:B_j\to B~|~j\in I)$ for the topology $J$, and for each
$j\in I$ morphisms $g_j:X_i\to X$, $g'_j:X_i\to X'$ in $\cC$
such that $Fg_j=Fg'_j=f_j$ and $(Lg_j)(t)=(Lg'_j)(t')$.
\end{lemma}
\begin{proof} As in the proof of proposition \ref{prop_filter-Deligne}
we shall say that a diagram of $F^{-1}B$ is a pair of finite sets
$(A,A')$ with $A\subset\Ob(F^{-1}B)$, $A'\subset\rMorph(F^{-1}B)$,
such that the source and target of every element of $A'$ lie in $A$.
We notice :

\begin{claim}\label{cl_cone-for-diagram}
Let $(A,A')$ be any diagram of $F^{-1}B$. There exist a covering
family $(f_j:B_j\to B~|~j\in I)$ and for every $j\in I$ an object
$Y_j$ of $F^{-1}B_j$ with a system of morphisms
$(g_{jX}:Y_j\to X~|~X\in A)$ such that the following holds. For
every $j\in I$, every $X,X',X''\in A$ and every $h:X\to X''$,
$h':X'\to X''$ in $A'$ we have $h\circ g_{jX}=h'\circ g_{jX'}$.
\end{claim}
\begin{pfclaim} Let $X_1,\dots,X_n$ be the elements of $A$.
Since $F$ is locally cofiltered, by a simple induction on $n$
we find a covering family $f'_\bullet:=(f'_{j'}:B_{j'}\to B~|~j'\in I')$
and morphisms $g'_{j'i}:Y'_{j'}\to X_i$ such that $Fg'_{j'i}=f'_{j'}$
for every $j'\in I'$ and $i=1,\dots,n$. Next, let
$T\subset A'\times A'$ be the subset of all pairs $(h,h')$
such that $h$ and $h'$ have the same target; we show that for
every $T'\subset T$, and every $j'\in I'$ there
exist a covering family $f''_{\bullet\bullet}:=
(f''_{j'\lambda}:B_{j'\lambda}\to B_{j'}~|~\lambda\in\Lambda_{j'})$
and morphisms $g''_{j'\lambda}:Y_{j'\lambda}\to Y'_{j'}$ with
$Fg''_{j'\lambda}=f''_{j'\lambda}$ and
$$
h\circ g'_{j'i}\circ g''_{j'\lambda}=h'\circ g'_{j'i'}\circ g''_{j'\lambda}
$$
for every $\lambda\in\Lambda_{j'}$, every $i,i',k=1,\dots,n$ and
every $(h:X_i\to X_k,h':X_{i'}\to X_k)$ in $T'$. The claim will
follow for $T'=T$, by letting $f_\bullet$ be the covering family
$(f'_j\circ f''_{j'\lambda}:B_{j'\lambda}\to B~|~
j'\in I,\ \lambda\in\Lambda_{j'})$, and letting $g_{\bullet\bullet}$
be the system of morphisms $g'_{j'i}\circ g''_{j'\lambda}$.

Now, if $T'=\emptyset$, there is nothing to show. Thus, let
$(h_0,h_1)\in T'$; by induction on the cardinality of $T'$,
we may assume that the assertion is known for the subset
$T'':=T'\setminus\{(h_0,h_1)\}$. After replacing $f'_\bullet$
by the covering family
$(f'_j\circ f''_{j'\lambda}~|~j'\in I,\ \lambda\in\Lambda_{j'})$
and the system $g'_{\bullet\bullet}$ with the system of morphisms
$g'_{j'i}\circ g''_{j'\lambda}$, we may then assume that
$h\circ g'_{j'i}=h'\circ g'_{j'i'}$ for every $i,i'=1,\dots,n$,
every $j\in I'$ and every $(h,h')\in T''$ such that the sources
of $h$ and $h'$ are $i$ and $i'$ respectively. We may assume
that the common target of $h_0$ and $h_1$ is $X_1$, and say
that their sources are $X_{i_0}$ and $X_{i_1}$ respectively.
Since $F$ is locally cofiltered, there exist a covering family
$f''_{\bullet\bullet}$ and a system of morphisms $g''_{\bullet\bullet}$
as in the foregoing, such that
$h_0\circ g'_{j'i_0}\circ g''_{j'\lambda}=
h_1\circ g'_{j'i_1}\circ g''_{j'\lambda}$ for every $j'\in I$ and
every $\lambda\in\Lambda_{j'}$. Clearly with these choices for
$f''_{\bullet\bullet}$ and $g''_{\bullet\bullet}$ the sought identities
hold whenever $(h,h')\in T'$, as required.
\end{pfclaim}

Now, the condition $[X,t]=[X',t']$ means that
for some $k\in\N$ there exist objects of $F^{-1}B$ 
$$
X_0,Y_0,X_1,Y_1,\dots,Y_k,X_{k+1}
$$
and sections $t_i\in LX_i$ for $i=0,\dots,k+1$, with $X_0=X$,
$X_{k+1}=X'$, $t_0=t$ and $t_{k+1}=t'$, as well as morphisms
$$
X_i\xleftarrow{q_i}Y_i\xrightarrow{q'_i}X_{i+1}
\qquad
\text{such that $(Lq_i)(t_i)=(Lq'_i)(t_{i+1})$ for $i=0,\dots,k$}.
$$
We apply claim \ref{cl_cone-for-diagram} to the diagram formed
by all the objects $X_i,Y_i$ and all the morphisms $q_i,q'_i$,
to find a covering family $(f_j:B_j\to B~|~j\in I)$ and morphisms
$h_j:Z_j\to Y_i$ with $Fh_j=f_j$, such that
$g_{ij}:=q_i\circ h_j=q'_i\circ h_j$ for every $j\in I$ and $i=0,\dots,k$.
It follows that $(Lg_{0j})(t_0)=(Lg_{1j})(t_1)=\cdots=(Lg_{k+1,j})(t_{k+1})$
for every $j\in I$. Then the assertion holds with $g_j:=g_{0j}$
and $g'_j:=g_{k+1,j}$ for every $j\in I$.
\end{proof}

\begin{proposition}\label{prop_sheaf-Fubini}
With the notation of \eqref{subsec_Fubini-sheaf}, we have :
\begin{enumerate}
\item
The functor $\int^\blambda$ commutes with all colimits.
\item
If $F$ is a locally cofiltered fibration, $\int^\blambda_a$ is exact.
\end{enumerate}
\end{proposition}
\begin{proof}(i): Since the colimits in $\cB^\wedge$ are computed
argumentwise (corollary \ref{cor_pre-misc}(ii)), we are reduced
to checking that for every $B\in\Ob(\cB)$ the functor
$$
\cC^\wedge\to\Set
\qquad
G\mapsto\colim_{(F^{-1}B)^o}G\circ\iota^o_B
$$
commutes with colimits. However, the functor
$\iota_B^\wedge:\cC^\wedge\to(F^{-1}B)^\wedge$ induced by the
inclusion functor $\iota_B:F^{-1}B\to\cC$ commutes with colimits
(again, by corollary \ref{cor_pre-misc}(ii)), so we come down to
checking that the functor $\Colim_{(F^{-1}B)^o}:(F^{-1}B)^\wedge\to\Set$
of remark \ref{rem_Lims-and-Colims}(ii) commutes with colimits,
which is clear, as the latter is a left adjoint. 

(ii): In light of (i), and since the functor $G\mapsto G^a$
commutes with colimits, it suffices to check that $\int^\blambda_a$
is left exact, and by proposition \ref{prop_commute-criteria}(i)
we are reduced to showing that $\int^\blambda_a$ commutes with
equalizers and finite non-empty products, and that it preserves
final objects.

Thus, let $E$ be the presheaf on $\cC$ such that $E(X)=\{\emptyset\}$
for every $X\in\Ob(\cC)$, and set $E':=\int^\blambda E$. For every
$B\in\Ob(\cB)$ and every $s\in E'^a(B)$ there is a covering family
$B_\bullet:=(B_i\to B~|~i\in I)$ in the topology $J$, and for every
$i\in I$ an object $X_i$ of $F^{-1}B_i$ such that $s$ is represented
by the system of section $([X_i,\emptyset])\in E'(B_i)$. Conversely,
every such system defines a section of $E'^a(B)$. Now, since $F$ is
locally cofiltered there exists such a covering family $B_\bullet$
of $B$ with the property that $\Ob(F^{-1}B_i)$ is not empty for
every $i\in I$; by picking arbitrary $X_i\in\Ob(F^{-1}B_i)$ for
each $i$ we obtain therefore a section of $E'^a(B)$; this shows
that $E'^a(B)\neq\emptyset$ for every $B\in\Ob(\cB)$. Next, let
$s:=([X_i,\emptyset]~|~i\in I)$ and $s':=([X'_i,\emptyset]~|~i\in I)$
be two sections of $E'^a(B)$, with $FX_i=FX'_i=B_i$ for every $i\in I$.
Since $F$ is locally cofiltered, we may find for every $i\in I$ a
covering family $(f_{i\lambda}:B_{i\lambda}\to B_i~|~\lambda\in\Lambda_i)$
and morphisms $g_{i\lambda}:Y_{i\lambda}\to X_i$,
$g'_{i\lambda}:Y_{i\lambda}\to X'_i$ in $\cC$ with
$Fg_{i\lambda}=Fg'_{i\lambda}=f_{i\lambda}$ for every $\lambda\in\Lambda_i$.
Then both $s$ and $s'$ are represented by the system
$([Y_{i\lambda},\emptyset]~|~i\in I,\ \lambda\in\Lambda_i)$. This
proves that $E'^a(B)$ contains a unique section for every
$B\in\Ob(\cB)$, {\em i.e.} $E'^a$ is the final object of $(\cB,J)^\sim$,
as required.

To check that $\int^\blambda_a$ commutes with finite products,
let $G_1,G_2$ be two presheaves on $\cC$, and set
$$
G:=G_1\times G_2
\qquad
H:=\int^\blambda G
\qquad
H_i:=\int^\blambda G_i
\qquad
\text{for $i=1,2$}.
$$
We need to check that the projections $G\to G_i$ induce an
isomorphism $\omega:H^a\isom H^a_1\times H^a_2$. To this aim,
let as well $B\in\Ob(\cB)$ and $s\in(H^a_1\times H^a_2)(B)$.
Hence $s$ is represented by the datum of a covering family
$(f_j:B_j\to B~|~j\in I)$ and sections $s_{ij}\in H_i(B_j)$
for $i=1,2$ and every $j\in I$, fulfilling the following
compatibility condition. For every $B'\in\Ob(\cB)$, every
$j,j'\in I$ and every pair of morphisms $g:B'\to B_j$,
$g':B'\to B_{j'}$ such that $f_j\circ g_j=f_{j'}\circ g'$, we
have $H_i(g)(s_{ij})=H_i(g')(s_{ij'})$ for $i=1,2$. By example
\ref{ex_complete-cats}(i), each $s_{ij}$ is the class
$[X_{ij},t_{ij}]$ of a pair consisting of an object $X_{ij}$ of
$F^{-1}B_j$ and a section $t_{ij}\in G_iX_{ij}$. Since
$F$ is locally cofiltered, we may then find for every $j\in I$
a covering family $(f_{j\lambda}:B_{j\lambda}\to B_j~|~\lambda\in\Lambda_j)$
and for every $j\in I$ and every $\lambda\in\Lambda_j$ an object
$Y_{j\lambda}$ of $F^{-1}B_{j\lambda}$ and morphisms
$g_{ij\lambda}:Y_{j\lambda}\to X_{ij}$ in $\cC$ with
$Fg_{ij\lambda}=f_{j\lambda}$  for $i=1,2$. With this notation,
set $t_{ij\lambda}:=(G_ig_{ij\lambda})(t_{ij})$; we have
$$
[Y_{j\lambda},t_{ij\lambda}]=H_i(f_{j\lambda})(s_{ij})
\qquad
\text{for $i=1,2$, every $j\in I$ and every $\lambda\in\Lambda_j$}.
$$
Thus, the class $\tau_{j\lambda}:=[Y_{j\lambda},(t_{1j\lambda},t_{2j\lambda})]$
is a section of $H(Y_{j\lambda})$ for every $j\in I$ and
$\lambda\in\Lambda_j$. Lastly, the family
$(f_j\circ f_{j\lambda}:B_{j\lambda}\to B~|~j\in I,\ \lambda\in\Lambda_j)$
covers $B$, and the system of sections $\tau_{\bullet\bullet}$ defines
a section of $H^a(B)$ whose image under $\omega_B$ agrees with
$s$. This shows that $\omega$ is an epimorphism. In order to check
that $\omega$ is a monomorphism, consider $t,t'\in H^a(B)$ whose
images agree in $(H_1^a\times H_2^a)(B)$; we may find a covering
family $f_\bullet:=(f_j:B_j\to B~|~j\in I)$ such that $t$ and $t'$
can be represented by compatible systems of sections
$t_\bullet:=([X_j,(t_{1j},t_{2j})]~|~j\in I)$,
$t'_\bullet:=([X'_j,(t'_{1j},t'_{2j})]~|~j\in I)$ with
$X_j,X'_j\in F^{-1}B_j$ and $t_{ij},t'_{ij}\in G_iX_j$ for $i=1,2$
and every $j\in I$. By assumption, the resulting systems
$([X_j,t_{ij}]~|~j\in I)$ and $([X'_j,t'_{ij}]~|~j\in I)$ represent
the same sections of $H^a_i(B)$ for $i=1,2$. This means that there
exists for every $j\in I$ a covering family
$(f_{j\lambda}:B_{j\lambda}\to B_j~|~\lambda\in\Lambda_j)$ such that
if we set $t_{ij\lambda}:=G_i(\blambda(f_{j\lambda},X_j))(t_{ij})$
and $t'_{ij\lambda}:=G_i(\blambda(f_{j\lambda},X'_j))(t'_{ij})$, we have
$$
[f^*_{j\lambda}X_j,t_{ij\lambda}]=[f^*_{j\lambda}X'_j,t'_{ij\lambda}]
\qquad
\text{in $H_i(B_{j\lambda})$ for $i=1,2$, every $j\in I$ and
every $\lambda\in\Lambda_j$}.
$$
After replacing $f_\bullet$ by the covering family
$(f_j\circ f_{j\lambda}:B_{j\lambda}\to B~|~j\in I,\ \lambda\in\Lambda_j)$
and the systems $t_\bullet$, $t'_\bullet$ by
$([f^*_{j\lambda}X_j,(t_{1j\lambda},t_{2j\lambda})]~|~j\in I,\ \lambda\in\Lambda)$
and
$([f^*_{j\lambda}X'_j,(t'_{1j\lambda},t'_{2j\lambda})]~|~j\in I,\ \lambda\in\Lambda)$,
we may therefore assume from start that
$$
[X_j,t_{ij}]=[X'_j,t'_{ij}]
\qquad
\text{in $H_i(B_j)$ for $i=1,2$ and every $j\in I$}.
$$
By lemma \ref{lem_locally-filt-equality}, we may then find for
$i=1,2$ and every $j\in I$ a covering family
$(f_{ij\lambda}:B_{ij\lambda}\to B_j~|~\lambda\in\Lambda_{ij})$ and morphisms
$g_{ij\lambda}:Y_{ij\lambda}\to X_j$, $g'_{ij\lambda}:Y_{ij\lambda}\to X'_j$ in
$\cC$ such that $Fg_{ij\lambda}=Fg'_{ij\lambda}=f_{ij\lambda}$ and
$t_{ij\lambda}:=(G_ig_{ij\lambda})(t_{ij})=(G_ig'_{ij\lambda})(t'_{ij})$
for every $\lambda\in\Lambda_{ij}$.

For $i=1,2$ and every $j\in I$, let $\cS_{ij}\subset\cB/B_j$ be the
sieve generated by the family $f_{ij\bullet}$, and pick a generating
family $(f'_\sigma:B_\sigma\to B_j~|~\sigma\in\Sigma_j)$ for the sieve
$\cS_{1j}\cap\cS_{2j}$, which covers $B_j$ as well in the topology
$J$. Set $\Sigma:=\bigcup_{j\in I}\{j\}\times\Sigma_j$. By definition,
this means that for $i=1,2$ and every $(j,\sigma)\in\Sigma$ there
exist $\lambda\in\Lambda_{1j}$, $\lambda'\in\Lambda_{2j}$ and a
commutative diagram in $\cB$ :
$$
\xymatrix{ B_\sigma \ar[r]^-{h_\sigma} \ar[d]_{h'_\sigma}
\ar[rd]^-{f'_\sigma} & B_{1j\lambda} \ar[d]^{f_{1j\lambda}} \\
B_{2j\lambda'} \ar[r]_-{f_{2j\lambda'}} & B_j.
}$$
With this notation, we set $Z_{1j\sigma}:=h^*_\sigma Y_{1j\lambda}$,
$Z_{2j\sigma}:=h'^*_\sigma Y_{2j\lambda'}$ and :
$$
s_{1j\sigma}:=G_1(\blambda(h_\sigma,Y_{1j\lambda}))(t_{1j\lambda})
\quad
s_{2j\sigma}:=G_2(\blambda(h'_\sigma,Y_{2j\lambda}))(t_{2j\lambda})
\qquad
\text{for every $(j,\sigma)\in\Sigma$}.
$$
Then, both compatible systems $([X_j,t_{1j}]~|~j\in I)$ and
$([X'_j,t'_{1j}]~|~j\in I)$ agree with the compatible system
$([Z_{1\sigma},s_{1\sigma}]~|~\sigma\in\Sigma)$ in $H^a_1(B)$, and
both compatible systems $([X_j,t_{2j}]~|~j\in I)$
and $([X'_j,t'_{2j}]~|~j\in I)$ agree with the compatible system
$([Z_{2\sigma},s_{2\sigma}]~|~\sigma\in\Sigma)$ in $H^a_2(B)$.

Since $F$ is locally cofiltered, we may find for every
$\sigma\in\Sigma$ a covering family
$(f_{\sigma\lambda}:B_{\sigma\lambda}\to B_\sigma~|~\lambda\in\Lambda_\sigma)$
and morphisms $g_{i\sigma\lambda}:Z_{\sigma\lambda}\to Z_{i\sigma}$, with
$Fg_{i\sigma\lambda}=f_{\sigma\lambda}$ for $i=1,2$.
Set $s_{i\sigma\lambda}:=(G_ig_{i\sigma\lambda})(s_{i\sigma})$ for $i=1,2$,
every $\sigma\in\Sigma$ and every $\lambda\in\Lambda_\sigma$. Hence
both $t$ and $t'$ are represented by the compatible system of
sections $([Z_{\sigma\lambda},(s_{1\sigma\lambda},s_{2\sigma\lambda})]~|~
\sigma\in\Sigma,\ \lambda\in\Lambda_\sigma)$, so $t=t'$.

Lastly, let $\phi_1,\phi_2:G_1\to G_2$ be two morphisms of presheaves
on $\cC$, and $E\subset G_1$ the equalizer of $\phi_1$ and $\phi_2$;
set also $H_i:=\int^\blambda G_i$, $\psi_i:=\int^\blambda\phi_i$ for
$i=1,2$, and denote by $E'$ the equalizer of
$\psi^a_1,\psi^a_2:H_1^a\to H_a^2$.
We need to check that the induced morphism $E^a\to E'$ is an
isomorphism. Thus, let $B\in\Ob(\cB)$ and $s\in E'(B)$; by remark
\ref{rem_rep-and-sheafify}(i), the set $E'(B)$ is the equalizer of
the induced maps $\psi^a_1(B),\psi^a_2(B):H^a_1(B)\to H^a_2(B)$, so
we may find a covering family $(f_j:B_j\to B~|~j\in I)$ such that
$s$ is represented by a compatible system
$t_\bullet:=([X_j,t_j]~|~j\in I)$ with $X_j\in\Ob(F^{-1}B_j)$ and
$t_j\in G_1(B_j)$ for every $j\in I$, such that the systems
$([X_j,(\phi_{1X_j})(t_j)]~|~j\in I)$ and
$([X_j,(\phi_{2X_j})(t_j)]~|~j\in I)$ represent the same section
of $H^a_2(B)$. This means that for every $j\in I$ there exists a
covering family
$f_\bullet:=(f_{j\lambda}:B_{j\lambda}\to B_j~|~\lambda\in\Lambda_j)$
such that if we set $t_{j\lambda}:=G_1(\blambda(f_{j\lambda},X_j))(t_j)$
we get the following identity in $H_2(B_{j\lambda})$ : 
$$
[f^*_{j\lambda}X_j,\phi_{1X_j}(t_{j\lambda})]=
[f^*_{j\lambda}X_j,\phi_{2X_j}(t_{j\lambda})]
\qquad
\text{for $i=1,2$, every $j\in I$ and every $\lambda\in\Lambda_j$}.
$$
After replacing $f_\bullet$ by the covering family
$(f_{j\lambda}\circ f_j:B_{j\lambda}\to B~|~j\in I,\ \lambda\in\Lambda_j)$
and $t_\bullet$ by the compatible system
$([f^*_{j\lambda}X_j,t_{j\lambda}]~|~j\in I,\ \lambda\in\Lambda_j)$
we may then assume that
$$
[X_j,(\phi_{1X_j})(t_j)]=[X_j,(\phi_{2X_j})(t_j)]
\qquad
\text{in $H_2(B_j)$ for every $j\in I$}.
$$
By lemma \ref{lem_locally-filt-equality} we may then find
for every $j\in I$ a covering family
$(f_{j\lambda}:B_{j\lambda}\to B_j~|~\lambda\in\Lambda_j)$ and
morphisms $g_{j\lambda},g'_{j\lambda}:Y_{j\lambda}\to X_j$ in $\cC$
such that $Fg_{j\lambda}=Fg'_{j\lambda}=f_{j\lambda}$ and
$$
G_2(g_{j\lambda})((\phi_{1X_j})(t_j))=
G_2(g'_{j\lambda})((\phi_{2X_j})(t_j))
\qquad
\text{for every $\lambda\in\Lambda_j$}.
$$
Set $\Sigma:=\bigcup_{j\in I}\{j\}\times\Lambda_j$; since $F$
is locally cofiltered, for every $\sigma\in\Sigma$ we may
then find a covering family
$(f'_{\sigma\lambda'}:B'_{\sigma\lambda'}\to B_\sigma~|~
\lambda'\in\Lambda'_\sigma)$ and a morphism
$h_{\sigma\lambda'}:Y_{\sigma\lambda'}\to Y_\sigma$ in $\cC$ such that
$Fh_{\sigma\lambda'}=f'_{\sigma\lambda'}$ and
$g_{\sigma\lambda'}:=g_\sigma\circ h_{\sigma\lambda'}=
g'_\sigma\circ h_{\sigma\lambda'}$ for every $\lambda'\in\Lambda'_\sigma$.
Set $t_{j\lambda\lambda'}:=G_1(g_{j\lambda\lambda'})(t_j)$ for every
$(j,\lambda)\in\Sigma$ and every $\lambda'\in\Lambda'_{(j,\lambda)}$.
The compatible system $([Y_{\sigma\lambda'},t_{\sigma\lambda'}]
~|~\sigma\in\Sigma,\ \lambda'\in\Lambda'_\sigma)$ represents
again the section $s$, and by construction we have
$\phi_{1Y_{\sigma\lambda'}}(t_{\sigma\lambda'})=
\phi_{2Y_{\sigma\lambda'}}(t_{\sigma\lambda'})$ for every $\sigma\in\Sigma$
and every $\lambda'\in\Lambda'_\sigma$. This proves that $s$ lies
in the image of the map $E^a(B)\to E'(B)$, so the natural morphism
$E^a\to E'$ is an epimorphism. To see that it is also a monomorphism,
consider two sections $s,s'$ of $E^a(B)$ whose images agree in
$E'(B)$; we may find a covering family $(f_j:B_j\to B~|~i\in I)$
and compatible systems $([X_j,t_j]~|~j\in I)$,
$([X'_j,t'_j]~|~j\in I)$ with $FX_j=FX'_j=B_j$ and $t_j\in G_1X_j$,
$t'_j\in G_1X'_j$ for every $j\in I$. Arguing as in the foregoing,
we may then assume that $[X_j,t_j]=[X'_j,t'_j]$ in $H_1(B_j)$ for
every $j\in I$. Then by lemma \ref{lem_locally-filt-equality}
there exist for every $j\in I$ a covering family
$(f_{j\lambda}:b_{j\lambda}\to B_j~|~j\in\Lambda_j)$ and morphims
$g_{j\lambda}:Y_{j\lambda}\to X_j$, $g'_{j\lambda}:Y_{j\lambda}\to X'_j$
in $\cC$ such that $Fg_{j\lambda}=Fg'_{j\lambda}=f_{j\lambda}$ and
$t_{j\lambda}:=(G_1g_{j\lambda})(t_j)=(G_1g'_{j\lambda})(t'_j)$ for every
$\lambda\in\Lambda_j$. Clearly
$\phi_{1Y_{j\lambda}}(t_{j\lambda})=\phi_{2Y_{j\lambda}}(t_{j\lambda})$ for
every $j\in I$ and every $\lambda\in\Lambda_j$; thus, both $s$
and $s'$ are represented by the compatible system
$([Y_{j\lambda},t_{j\lambda}]~|~j\in I,\ \lambda\in\Lambda_j)$,
and the proof is concluded.
\end{proof}

\sset\subsubsection{}\label{subsec_Fubini-is-g-shrieck}
Let $C:=(\cC,J)$ and $C':=(\cC',J')$ be two sites, $g:\cC\to\cC'$
a given functor, and $\sV$ a universe such that $C$ and $C'$ are
$\sV$-small. We shall apply proposition \ref{prop_sheaf-Fubini}
to the split fibration $\ss:\cC'/g\cC\to\cC'$ associated with $g$,
with its canonical cleavage $\blambda$ (see example
\ref{ex_split-fibration}). Let also $\st:\cC'/g\cC\to\cC$ be
the target functor (see \eqref{subsec_target-fctr}); with this
notation, notice that
$$
g_{\sV!}=\int^\blambda\circ\ \ \st^\wedge.
$$
In light of proposition \ref{prop_sheaf-Fubini}, it follows
already that if $\ss$ is locally cofiltered, then the functors
$g^a_{\sV!}$ and $\tilde g{}^*_\sV:C^\sim_\sV\to C'^\sim_\sV$ are
exact. In fact, we have :

\begin{proposition}\label{prop_morph-of-sites}
With the notation of \eqref{subsec_Fubini-is-g-shrieck},
the following conditions are equivalent :
\begin{enumerate}
\alphaenu
\item
$g$ is a morphism of sites $C'\to C$.
\item
The fibration $\ss$ is locally cofiltered, and for every
covering family $(f_i:X_i\to X~|~i\in I)$ for the topology $J$,
the family $(g(f_i):gX_i\to gX~|~i\in I)$ covers $gX$ relative
to $J'$.
\end{enumerate}
\end{proposition}
\begin{proof} (a)$\Rightarrow$(b) : Since $g$ is continuous,
we know already that it transforms covering families for $J$
into covering families for $J'$ (lemma \ref{lem_crit-continuity}).
Notice as well that $g^a_{\sV!}=\tilde g{}^*_\sV\circ(-)^a$
(lemma \ref{lem_cont-funct-site}(ii)); since $\tilde g{}^*_\sV$
is exact, the same holds then for $g^a_{\sV!}$.

Let us check next that condition (a) of definition
\ref{def_locally-cofiltered} holds for $\ss$. To this aim,
consider the final object $E$ of $C^\sim$ such that
$E(X)=\{\emptyset\}$ for every $X\in\Ob(\cC)$. By assumption,
$\tilde g{}^*_\sV(E)$ is the final object of $C'^\sim_\sV$,
{\em i.e.} $(\tilde g{}^*_\sV E)(X')$ is a set with one
element for every $X'\in\Ob(\cC')$. However, every section
of $(\tilde g{}^*_\sV E)(X')$ is represented by the datum
of a covering family $(X'_j\to X'~|~j\in I)$ in the
topology $J'$, and a system of sections
$([X'_j\to gX_j,\emptyset]~|~j\in I)$; especially,
$\Ob(\ss^{-1}X'_j)\neq\emptyset$ for every $j\in I$, as
required.

Next, let $l_i:X'\to gX_i$ for $i=1,2$ be two objects of
$\ss^{-1}X'$; in order to check condition (b) of definition
\ref{def_locally-cofiltered} we need to exhibit a covering
family $(f_j:X'_j\to X'~|~j\in I)$ and for every $j\in I$
an object $h_j:X'_j\to gY_j$ of $\ss^{-1}X'_j$ and morphisms
$k_{ij}:Y_j\to X_i$ in $\cC$ that make commute the diagram
\set\begin{equation}\label{eq_condition-b-cofilt}
{\diagram
X' \ar[rd]_-{l_i} & \ar[l]_-{f_j} X'_j \ar[r]^-{h_j}
& gY_j \ar[ld]^{g(k_{ij})} \\
& gX_i
\enddiagram}
\qquad
\text{for $i=1,2$}.
\end{equation}
To this aim, set $H_i:=g_{\sV!}(h_{X_i})=h_{gX_i}$ for $i=1,2$,
and $H:=g_{\sV!}(h_{X_1}\times h_{X_2})$. The pair $(l_1,l_2)$
yields a section $l_\bullet$ of $(H_1\times H_2)(X')$; on
the other hand, since $g^a_{\sV!}$ is exact, the natural
morphism $\omega:H\to H_1\times H_2$ induces an isomorphism
on associated sheaves in $C'^\sim$. Especially, there exists
a covering family $(f_j~|~j\in I)$ as sought, such that
for every $j\in I$ the section
$(H_1\times H_2)(f_j)(l_\bullet)=(l_1\circ f_j,l_2\circ f_j)$
agrees with the image under $\omega$ of a section
$[h_j:X'_j\to gY_j,(k_{1j},k_{2j})]$ of $H(X'_j)$. Unwinding
the definition, we obtain precisely a diagram
\eqref{eq_condition-b-cofilt}.

Lastly, let $l_i:X'\to gX_i$ for $i=1,2$ be two objects
of $\ss^{-1}X'$, and $t_1,t_2:X_1\to X_2$ two morphisms in
$\cC$ such that $g(t_1)\circ l_1=l_2=g(t_2)\circ l_1$; in
order to verify condition (c) of definition
\ref{def_locally-cofiltered} we need to exhibit a covering
family $(f_j:X'_j\to X'~|~j\in I)$ and for every $j\in I$
an object $h_j:X'_j\to gY_j$ of $\ss^{-1}X'_j$ with a
morphism $k_j:Y_j\to X_1$ such that
\set\begin{equation}\label{eq_cond-c-cofilt}
l_1\circ f_j=g(k_j)\circ h_j
\qquad\text{and}\qquad
t_1\circ k_j=t_2\circ k_j.
\end{equation}
To this aim, let $t_{i*}:h_{X_1}\to h_{X_2}$ be the morphism
of presheaves on $\cC$ induced by $t_i$, for $i=1,2$; we
set $H_i:=g_{\sV!}(h_{X_i})=h_{gX_i}$ for $i=1,2$, and denote
by $E$ the equalizer of $t_{1*}$ and $t_{2*}$, and by $E'$
the equalizer of $g_{\sV!}(t_{1*}),g_{\sV!}(t_{2*}):H_1\to H_2$.
Since $g^a_{\sV!}$ is exact, the natural morphism
$\omega:g_{\sV!}E\to E'$ induces an isomorphism on
associated sheaves in $C^\sim$.
On the other hand, $l_1$ defines a section of $E'(X')$;
it follows that there exists a covering family
$(f_j:X'_j\to X'~|~j\in I)$ such that for every $j\in I$
the section $(E'f_j)(l_1)$ agrees with the image under
$\omega$ of a section $[h_j:X'_j\to gY_j,k_j]$ of
$(g_{\sV!}E')(X'_j)$. Unwinding the definitions, we obtain
the identities \eqref{eq_cond-c-cofilt}.

(b)$\Rightarrow$(a): We have already noticed that if $\ss$
is locally cofiltered, the functor $\tilde g{}^*_\sV$ is exact.
In order to check that $g$ is continuous, let us consider
any sheaf $F$ on $\cC'$, any $X\in\Ob(\cC)$, and any
covering subobject $R\subset h_X$ for the topology $J$.
We need to show that the natural map
$$
F(gX)\simeq\Hom_{\cC^\wedge_\sV}(h_X,g^*_\sV F)\to
\Hom_{\cC^\wedge_\sV}(R,g^*_\sV F)
$$
is bijective. By adjunction, the latter is naturally
identified with the induced map
$$
\Hom_{\cC'^\wedge_\sV}(g_{\sV!}(h_X),F)\to
\Hom_{\cC'^\wedge_\sV}(g_{\sV!}(R),F)
$$
so we come down to checking that the morphism
$\phi:(g_{\sV!}R)^a\to(g_{\sV!}h_X)^a=h^a_{gX}$ induced by the
inclusion $R\to h_X$ is an isomorphism. However, we know
already that the functor $g^a_{\sV!}$ is exact, so $\phi$
is a monomorphism. Next, let $(f_j:X_j\to X~|~j\in I)$ be
a $\sV$-small family of generators for the covering sieve
of $X$ corresponding to $R$; there follows a covering
morphism
$$
S:=\coprod_{j\in I}h_{X_j}\to R\to h_X
\qquad
\text{in $\cC^\wedge$}
$$
whose image under $g_{\sV!}$ is still a covering morphism :
indeed, $g_{\sV!}$ commutes with coproducts since it is a
left adjoint, so $g_{\sV!}(S)$ is the coproduct of the family
$(h_{gX_j}~|~j\in I)$, and by assumption $g$ transforms covering
families for the topology $J$ into covering families for the
topology $J'$. Thus, the morphism $g_{\sV!}R\to g_{\sV!}h_X$
induces an epimorphism as well on associated sheaves, whence
the contention.
\end{proof}

\begin{example}\label{ex_simple-case}
In the situation of \eqref{subsec_Fubini-is-g-shrieck},
suppose that $C$ is a lex-site and $g$ is left exact.
Then $g$ is a morphism of sites if and only if for every
covering family $(f_i:X_i\to X~|~i\in I)$ for the topology $J$,
the family $(g(f_i):gX_i\to gX~|~i\in I)$ covers $gX$ relative
to $J'$. Indeed, under this assumption, the category
$X/g\cC$ is cofiltered for every $X\in\Ob(\cC')$ (example
\ref{ex_cofiltered-comma}(i)), so the source fibration
$\ss:\cC'/g\cC\to\cC'$ is trivially locally cofiltered, and
the assertion follows from proposition \ref{prop_morph-of-sites}.
\end{example}

\begin{theorem}\label{th_gener-Beilinson}
Let $(\cC,J)$ be a $\sU$-site, $\cC'$ a category, $g:\cC'\to\cC$
a functor, and endow $\cC'$ with the topology $J'$ induced by
$g$. Suppose that the following conditions hold :
\begin{enumerate}
\alphaenu
\item
For every $X\in\Ob(\cC)$ there exists a family
$(gY_i\to X~|~i\in I)$ of objects of $g\cC'/X$ that covers $X$
in the topology $J$.
\item
For every $Y,Y'\in\Ob(\cC')$ and every morphism $\phi:gY\to gY'$
in $\cC$ there exist :
\begin{itemize}
\item[(b.i)]
a family $(\psi_i:Z_i\to Y~|~i\in I)$ of objects of\/ $\cC'/Y$ such
that the family $(g(\psi_i):gZ_i\to gY~|~i\in I)$ covers $gY$ in
the topology $J$
\item[(b.ii)]
for each $i\in I$ a morphism $\nu_i:Z_i\to Y'$ such that
$\phi\circ g(\psi_i)=g(\nu_i)$.
\end{itemize}
\item
For every pair of morphisms $\phi,\phi':Y\to Y'$ in $\cC'$ such that
$g(\phi)=g(\phi')$ there exists a family $(\psi_i:Z_i\to Y~|~i\in I)$
of objects of\/ $\cC'/Y$ such that the family
$(g(\psi_i):gZ_i\to gY~|~i\in I)$ covers $gY$ in the topology $J$,
and $\phi\circ\psi_i=\phi'\circ\psi_i$\ \ for every $i\in I$.
\end{enumerate}
Then we have :
\begin{enumerate}
\item
$g$ is cocontinuous for the topologies $J,J'$ and is a
morphism of sites $(\cC,J)\to(\cC',J')$.
\item
$g$ induces an equivalence $\tilde g_*:(\cC,J)^\sim\isom(\cC',J')^\sim$.
\item
A family $(\psi_i:Y_i\to Y~|~i\in I)$ of morphisms of\/ $\cC'$ generates
a covering sieve of $J'$ if and only if $(g(\psi_i):gY_i\to gY~|~i\in I)$
generates a covering sieve of $J$.
\end{enumerate}
\end{theorem}
\begin{proof} Suppose first that both $\cC$ and $\cC'$ are
small. Let $\ss:\cC/g\cC'\to\cC$ be the fibration associated
with $g$, as in \eqref{subsec_Fubini-is-g-shrieck}; we remark :

\begin{claim}\label{cl_abc-imply-s-loc-cofilt}
The fibration $\ss$ is locally cofiltered.
\end{claim}
\begin{pfclaim} To check condition (a) of definition
\ref{def_locally-cofiltered}, let $X\in\Ob(\cC)$; we need
to find a covering family $(X_j\to X~|~j\in I)$ for the
topology $J$, such that $\Ob(\ss^{-1}X_j)\neq\emptyset$
for every $j\in I$. However, by condition (a) of the theorem
we have a covering family $(gY_j\to X~|~j\in I)$ for $X$; we
may then choose $X_j:=gY_j$, since then
$\one_{gY_j}\in\Ob(\ss^{-1}X_j)$ for every $j\in I$.

Next, we check condition (b) of definition
\ref{def_locally-cofiltered} : consider any $X\in\Ob(\cC)$
and morphisms $\phi_i:X\to gY_i$ in $\cC'$ for $i=1,2$; we
need to find a covering family $(f_j:X_j\to X~|~j\in I)$
such that for every $j\in I$ there exist a morphism
$\psi_j:X_j\to gZ_j$ in $\cC$ and morphisms $\nu_{ij}:Z_j\to Y_i$
with $g(\nu_{ij})\circ\psi_j=\phi_i\circ f_j$ for $i=1,2$.
To this aim, we use first condition (a) of the theorem to
find a covering family $f'_\bullet:=(f'_j:gY'_j\to X~|~j\in I')$;
then, by condition (b) of the theorem we find for every $j\in I'$
a family $(f'_{j\lambda}:Y'_{j\lambda}\to Y'_j~|~\lambda\in\Lambda_j)$
such that $(g(f'_{j\lambda})~|~\lambda\in\Lambda_j)$ covers $gY'_j$,
and such that for every $\lambda\in\Lambda_j$ there exists a
morphism $h_{j\lambda}:Y'_{j\lambda}\to Y_1$ with
$g(h_{j\lambda})=\phi_1\circ f'_j\circ g(f'_{j\lambda})$. We may
then replace $f'_\bullet$ by the family
$(f'_j\circ g(f'_{j\lambda}):gY'_{j\lambda}\to X~|~
j\in I',\ \lambda\in\Lambda_j)$, and assume from start that for
every $j\in I'$ there exists a morphism $h_j:Y'_j\to Y_1$ such
that $g(h_j)=\phi_1\circ f'_j$. Next, we apply again condition
(b) of the theorem to find for every $j\in I'$ a family
$(f''_{j\lambda'}:Y''_{j\lambda'}\to Y'_j~|~\lambda'\in\Lambda'_j)$
such that $(g(f''_{j\lambda'})~|~\lambda'\in\Lambda'_j)$ covers
$gY'_j$, and such that for every $\lambda'\in\Lambda'_j$ there
exists a morphism $h'_{j\lambda'}:Y''_{j\lambda'}\to Y_2$ with
$g(h'_{j\lambda'})=\phi_2\circ f'_j\circ g(f''_{j\lambda'})$.
We may then further replace $f'_\bullet$ by the family
$(f'_j\circ g(f''_{j\lambda'}):gY''_{j\lambda'}\to X~|~j\in I',\
\lambda'\in\Lambda'_j)$, and the system $(h_j~|~j\in I')$
by the system $(h_j\circ f''_{j\lambda'}:Y''_{j\lambda'}\to Y_1~|~
j\in I',\ \lambda'\in\Lambda'_j)$, and assume from start that
there exists as well for every $j\in I'$ a morphism
$h'_j:Y'_j\to Y_2$ such that $g(h'_j)=\phi_2\circ f'_j$.
Then we set $X_j:=gY'_j$, $Z_j:=Y'_j$, $\psi_j:=\one_{X_j}$,
$\nu_{1j}:=h_j$ and $\nu_{2j}:=h'_j$ for every $j\in I'$;
clearly the resulting family $f_\bullet$ and systems of
morphisms $\psi_\bullet$, $\nu_{\bullet\bullet}$ will do.

Lastly, we check condition (c) of definition
\ref{def_locally-cofiltered} : let $X\in\Ob(\cC)$ and
consider morphisms $\phi:X\to gY_1$ in $\cC$ and
$\psi,\psi':Y_1\to Y_2$ in $\cC'$ such that
$g(\psi)\circ\phi=g(\psi')\circ\phi$; we need to find
a covering family $(f_j:X_j\to X~|~j\in I)$ and for every
$j\in I$ morphisms $h_j:X_j\to gZ_j$ in $\cC$ and
$\nu_j:Z_j\to Y_1$ in $\cC'$ such that
$\phi\circ f_j=g(\nu_j)\circ h_j$ and
$\psi\circ\nu_j=\psi'\circ\nu_j$. To this aim, we use
conditions (a) and (b) of the theorem to get first a
covering family $f'_\bullet:=(f'_j:gY_j\to X~|~j\in I)$,
and then for every $j\in I$ a family
$(f''_{j\lambda}:Y_{j\lambda}\to Y_j~|~\lambda\in\Lambda_j)$
such that $(g(f''_{j\lambda})~|~\lambda\in\Lambda_j)$ covers
$gY_j$, and such that for every $\lambda\in\Lambda_j$ there
exists a morphism $\nu'_{j\lambda}:Y_{j\lambda}\to Y_1$ with
$g(\nu'_{j\lambda})=\phi\circ f'_j\circ g(f''_{j\lambda})$.
After replacing $f'_\bullet$ by the family
$(f'_j\circ g(f''_{j\lambda})~|~j\in I,\ \lambda\in\Lambda_j)$,
we may assume that for every $j\in I$ there exists a morphism
$\nu'_j:Y_j\to Y_1$ such that $g(\nu'_j)=\phi\circ f'_j$.
Especially, notice that
$g(\psi\circ\nu'_j)=g(\psi'\circ\nu'_j)$ for every $j\in I$.
By condition (c) of the theorem, there exists therefore
for every $j\in I$ a family
$(\nu'_{j\lambda'}:Y'_{j\lambda'}\to Y_j~|~\lambda'\in\Lambda'_j)$
such that $(g(\nu'_{j\lambda'})~|~\lambda'\in\Lambda'_j)$
covers $gY_j$ and $\psi\circ\nu'_j\circ\nu'_{j\lambda'}=
\psi'\circ\nu'_j\circ\nu_{j\lambda}$ for every
$\lambda'\in\Lambda_j$. Set
$I:=\bigcup_{j\in I}\{j\}\times\Lambda_j$; for every
$(j,\lambda')\in I$ we let $X_{j\lambda'}:=gY_{j\lambda'}$,
$f_{j\lambda'}:=f'_j\circ g(\nu'_{j\lambda'})$,
$h_{j\lambda'}:=\one_{X_{j\lambda'}}$ and
$\nu_{j\lambda'}:=\nu'_j\circ\nu'_{j\lambda'}$. Clearly the
resulting covering family $f_{\bullet\bullet}$ and systems
of morphisms $h_{\bullet\bullet}$, $\nu_{\bullet\bullet}$ will do.
\end{pfclaim}

(iii) follows from claim \ref{cl_abc-imply-s-loc-cofilt}, lemmata
\ref{lem_induced-top}(ii) and \ref{lem_crit-continuity}, and
proposition \ref{prop_morph-of-sites}.

(i): The continuity of $g$ holds by definition of $J'$;
then $g$ is a morphism of sites, by claim
\ref{cl_abc-imply-s-loc-cofilt} and proposition
\ref{prop_morph-of-sites}. To check the cocontinuity, let
$Y\in\Ob(\cC')$ and $f_\bullet:=(f_j:X_j\to gY~|~j\in I)$ a
family of morphisms in $\cC$ generating a covering sieve
$\cS$ of $gY$; we need to show that $g^{-1}_{|Y}\cS$ is a
covering sieve of $Y$ for the topology $J'$. However, by
(a) we may find for every $j\in I$ a family
$(f_{j\lambda}:gY_{j\lambda}\to X_j~|~\lambda\in\Lambda_j)$
covering $X_j$ in the topology $J$; let $\cS'\subset\cS$ be
the sieve generated by the system
$(f'_{j\lambda}:=f_j\circ f_{j\lambda}:gY_{j\lambda}\to gY~|~j\in I,\
\lambda\in\Lambda_j)$. Clearly $\cS'$ covers $gY$ in the topology
$J$, and by remark \ref{rem_topology}(ii) it suffices to check
that $g^{-1}_{|Y}\cS'$ is a covering sieve for $Y$, so we may
replace $f_\bullet$ by the family $f'_{\bullet\bullet}$, and assume
from start that for every $j\in I$ we have $X_j=gY_j$ for some
$Y_j\in\Ob(\cC')$. By condition (b), we may then find for every
$j\in I$ a family
$(h_{j\lambda}:Y_{j\lambda}\to Y_j~|~\lambda\in\Lambda_j)$ of morphisms
in $\cC'$ such that the family $(g(h_{j\lambda})~|~\lambda\in\Lambda_j)$
covers $gY_j$ in the topology $J$, and for every $j\in I$ and
$\lambda\in\Lambda_j$ a morphism $k_{j\lambda}:Y_{j\lambda}\to Y$ in
$\cC'$ such that $f''_{j\lambda}:=f_j\circ g(h_{j\lambda})=g(k_{j\lambda})$.
Let $\cS''\subset\cS$ be the sieve generated by the system
$(f''_{j\lambda}:gY_{j\lambda}\to gY~|~j\in I,\ \lambda\in\Lambda_j)$.
Then $\cS''$ covers $gY$ in the topology $J$, and it suffices
to check that $\cT:=g^{-1}_{|Y}\cS''$ covers $Y$ in the topology $J'$.
However, $\cT$ contains the family
$(k_{j\lambda}~|~j\in I,\ \lambda\in\Lambda_j)$, which covers $Y$,
by virtue of (iii). The assertion follows.

\begin{claim}\label{cl_just-for-units}
(i)\ \
Let $\phi:G\to G'$ be a morphism of presheaves on $\cC$
such that $g^\wedge(\phi):g^\wedge G\to g^\wedge G'$ is an
isomorphism in $\cC'^\wedge$. Then $\phi$ is a bicovering
morphism.

(ii)\ \
Let $(\eta',\eps')$ be a unit and counit for the adjoint
pair $(\tilde g{}^*,\tilde g_*)$. In order to prove
assertion (ii) of the theorem, it suffices to check that
$\eta':\one_{(\cC',J')^\sim}\Rightarrow\tilde g_*\tilde g{}^*$ is
an isomorphism of functors.
\end{claim}
\begin{pfclaim}(i): Indeed, let $X\in\Ob(\cC)$ and
$s\in G'X$; by condition (a) there exists a covering
family $(f_i:gY_i\to X~|~i\in I)$, and by assumption
$(G'f_i)(s)\in\Img(\phi_{gY_i})$ for every $i\in I$, so
$\phi$ is a covering morphism, by remark \ref{rem_iprippi}(ii).
Likewise, let $s,s'\in GX$ such that $t:=\phi_X(s)=\phi_X(s')$;
it follows that
$\phi_{gY_i}((Gf_i)(s))=\phi_{gY_i}((Gf_i)(s'))=(G'f_i)(t)$,
whence $(Gf_i)(s)=(Gf_i)(s')$ for every $i\in I$. In light
of remark \ref{rem_iprippi}(iii), the assertion follows.

(ii): Indeed, if $\eta'$ is an isomorphism, the triangular
identities of \eqref{subsec_adj-pair} show that
$g^\wedge(\eps'_G)$ will also be an isomorphism for every
sheaf $G$ on $(\cC,J)$, hence $\eps'_G$ will be an
isomorphism, by (i), and to conclude it will suffice to
invoke proposition \ref{prop_fullfaith-adjts}(i,iii).
\end{pfclaim}

(ii): From (i) and lemma \ref{lem_improve}(ii) we get an
isomorphism of functors $\cC^\wedge\to(\cC',J')^\sim$
$$
\omega:(-)^a\circ g^\wedge\isom\tilde g_*\circ(-)^a.
$$
Let now $\eta:\one_{\cC'^\wedge}\Rightarrow g^\wedge\circ g_!$
and $\eps:g_!\circ g^\wedge\Rightarrow\one_{\cC^\wedge}$ be the
unit and counit of the natural adjunction for the pair of
functors $(g_!,g^\wedge)$, as in remark \ref{rem_was-cofinal}(iii).
We consider the natural transformations of functors on
$(\cC',J')^\sim$ and respectively on $(\cC,J)^\sim$
$$
\eta'_F:F\xrightarrow{\ (\eta_F)^a\ }
(g^\wedge g_!F)^a\xrightarrow{\ \omega_{g_!F}\ }
(g^\wedge\circ g^a_!)(F)
\qquad
\eps'_G:(g^a_!\circ g^\wedge)(G)\xrightarrow{\ (\eps_F)^a\ }G.
$$
A little diagram chase shows that the pair $(\eta',\eps')$
fulfills the triangular identities of \eqref{subsec_adj-pair},
so these are the unit and counit of an adjunction for the
pair $(\tilde g{}^*,\tilde g_*)$. By claim
\ref{cl_just-for-units}(ii), we are thus reduced to showing
that $\eta_F$ is a bicovering morphisms, for every presheaf
$F$ on $\cC'$.

We check first that $\eta_F$ is a covering morphism : thus,
let $X\in\Ob(\cC'^\wedge)$, and $s\in(g^\wedge g_!F)(X)$; then
$s$ is a class $[\phi:gX\to gY,\sigma]$, where $\phi$ is a
morphism in $\cC$, and $\sigma\in FY$. By condition (b) there
exists a family $(f_j:X_j\to X~|~j\in I)$ of morphisms of $\cC'$
such that $(\phi(f_j)~|~j\in I)$ is a covering family for the
topology $J$, and such that for every $j\in I$ there exists
a morphism $h_j:X_i\to Y$ with $g(h_j)=\phi\circ g(f_j)$.
Set $\sigma_j:=(Fh_j)(\sigma)\in FX_i$ for every $j\in I$;
then $\eta_{F,X_j}(\sigma_j)$ is the class $[\one_{gX_j},\sigma_j]$
in $(g^\wedge g_!F)(X_j)$ for every $j\in I$. On the other hand,
$(g^\wedge g_!F)(f_j)(s)=[\phi\circ g(f_j):gX_j\to gY,\sigma]$
for every $j\in I$. According to remark \ref{rem_iprippi}(ii),
it suffices then to show that
$[\one_{gX_j},\sigma_j]=[\phi\circ g(f_j),\sigma]$ in
$(g^\wedge g_!F)(X_j)$, for every $j\in I$. But $h_j$ yields a
morphism $gX_j/h_j:(X_j,\one_{gX_j})\to(Y,\phi\circ g(f_j))$ in
$gX_j/g\cC'$ whence the assertion.

Next, let $t,t'\in FX$ be two sections such that
$\eta_{F,X}(t)=\eta_{F,X}(t')$, {\em i.e.} $[\one_{gX},t]=[\one_{gX},t']$
in $(g_!F)(gX)$. According to lemma \ref{lem_locally-filt-equality},
there exist a covering family $(f_j:gX_j\to gX~|~j\in I)$
for the topology $J$, and morphisms $\phi_j:gX_j\to gY_j$ in
$\cC$ and $h_jh'_j:Y_j\to X$ in $\cC'$ such that
$$
g(h_j)\circ\phi_j=g(h'_j)\circ\phi_j=f_j
\qquad\text{and}\qquad
(Fh_j)(t)=(Fh'_j)(t')
\qquad
\text{for each $j\in I$}.
$$
Then, by condition (b), we may find for every $j\in I$ a
family $(f'_{j\lambda}:X_{j\lambda}\to X_j~|~\lambda\in\Lambda_j)$ such
that $(g(f'_{j\lambda})~|~\lambda\in\Lambda_j)$ is a covering
family for the topology $J$, and such that for every
$\lambda\in\Lambda_j$ there exists a morphism
$\nu_{j\lambda}:X_{j\lambda}\to Y_j$ with
$g(\nu_{j\lambda})=\phi_j\circ g(f'_{j\lambda})$. After replacing
$(f_j~|~j\in\Lambda)$ by the covering family
$(f_j\circ g(f'_{j\lambda}):X_{j\lambda}\to X~|~j\in I,\ \lambda\in\Lambda)$,
and $(\phi_j~|~j\in J)$ by the system of morphisms
$(\phi_j\circ g(f'_{j\lambda})~|~j\in I,\ \lambda\in\Lambda_j)$,
we may then assume that for every $j\in I$ there exists a
morphism $\nu_j:X_j\to Y_j$ such that $\phi_j=g(\nu_j)$, in
which case we have $g(h_j\circ\nu_j)=g(h'_j\circ\nu_j)=f_j$ for
every $j\in I$. Then, by condition (c) there exists for every
$j\in I$ a family
$(f''_{j\lambda'}:X'_{j\lambda'}\to X_j~|~\lambda'\in\Lambda'_j)$
of morphisms in $\cC'$ such that
$(g(f''_{j\lambda'})~|~\lambda'\in\Lambda'_j)$ is a covering family
for the topology $J$, and such that
$h_j\circ\nu_j\circ f''_{j\lambda'}=h'_j\circ\nu_j\circ f''_j$.
After replacing $(f_j~|~j\in\Lambda)$ by the covering family
$(f_j\circ g(f''_{j\lambda'}):X_{j\lambda'}\to X~|~
j\in I,\ \lambda'\in\Lambda'_j)$ and $(\nu_j~|~j\in I)$ by the
system $(\nu_j\circ f''_{j\lambda'}~|~j\in I,\ \lambda'\in\Lambda'_j)$,
we may then assume that $\mu_j:=h_j\circ\nu_j=h'_j\circ\nu_j$ for
every $j\in I$. It follows that $(F\mu_j)(t)=(F\mu_j)(t')$ for
every $j\in I$; lastly, the family $(\mu_j:X_j\to X~|~j\in I)$
covers $X$ in the topology $J'$, by virtue of (iii). In view
of remark \ref{rem_iprippi}(iii), we conclude that $\eta_F$
is a bicovering morphism, as stated.

This completes the proof of the theorem in case $\cC$ and
$\cC'$ are small. Lastly, consider the case where $(\cC,J)$
is a $\sU$-site and $\cC'$ is a general category, and pick
a universe $\sV$ such that $\cC$ and $\cC'$ are $\sV$-small.
The theorem then applies to the functor $\tilde g_{\sV*}$,
so the latter is an equivalence, and we also get assertions
(i) and (iii), which are independent of the universe $\sU$.
It remains therefore only to check that $\tilde g_{\sU*}$ is
an equivalence, and taking into account the commutative
diagram
$$
\xymatrix{ (\cC,J)^\sim_\sU \ar[r]^-{\tilde g_{\sU*}} \ar[d] &
(\cC',J')^\sim_\sU \ar[d] \\
(\cC,J)^\sim_\sV \ar[r]^-{\tilde g_{\sV*}} & (\cC',J')^\sim_\sV 
}$$
(whose vertical arrows are the fully faithful inclusion
functors), we already see that $\tilde g_{\sU*}$ is fully
faithful. We notice :

\begin{claim}\label{cl_reduce-size}
For every covering family $(f_i:X_i\to X~|~i\in I)$ in the
topology $J$ there exists a small subset $I'\subset I$ such
that $(f_i:X_i\to X~|~i\in I')$ is still a covering family.
\end{claim}
\begin{pfclaim} Pick a small topologically generating family
$G\subset\Ob(\cC)$ for the site $(\cC,J)$; then for every
$i\in I$ there exists a covering family
$(f_{i\lambda}:X_{i\lambda}\to X_i~|~\lambda\in\Lambda_i)$ with
$X_{i\lambda}\in G$ for every $\lambda\in\Lambda_i$. Set
$\Lambda:=\bigcup_{i\in I}\{i\}\times\Lambda_i$; the family
$(f'_{i\lambda}:=f_i\circ f_{i\lambda}:X_{i\lambda}\to X~|~
(i,\lambda)\in\Lambda)$ covers $X$ in the topology $J$. Since
$\Hom_\cC(X_{i\lambda},X)$ is a small set for every
$(i,\lambda)\in\Lambda$, we may then find a small subset
$\Lambda'\subset\Lambda$ such that
$(f'_{i\lambda}~|~(i,\lambda)\in\Lambda')$ is still covering;
then we can take for $I'\subset I$ the image of $\Lambda'$
under the natural projection $\Lambda\to I$.
\end{pfclaim}

Now, let $F$ be any $\sU$-sheaf on $(\cC',J')$; we know already
that there exists a $\sV$-sheaf $G$ on $(\cC,J)$ such that
$\tilde g_{*\sV}G$ is isomorphic to $F$, and we need to check
that $GX$ is essentially small for every $X\in\Ob(\cC)$.
But by condition (a) we may find a covering family
$(gY_i\to X~|~i\in I)$, and by claim \ref{cl_reduce-size}
we may assume that $I$ is small. Then the induced map
$GX\to\prod_{i\in I}G(gY_i)$ is injective, and by assumption
$G(gY_i)=FY_i$ is small for every $i\in I$, whence the contention.
\end{proof}

\begin{remark}\label{rem_gener-Beilinson}
(i)\ \
Notice that if $g$ is full (resp. faithful), then condition (b)
(resp. (c)) of theorem \ref{th_gener-Beilinson} is trivially
satisfied.

(ii)\ \
It follows easily that theorem \ref{th_gener-Beilinson} generalizes
\cite[\S2.1, Prop.]{Bei}. One can also show that conditions (a),
(b) and (c) of the theorem are implied by conditions $(\cL0)$, $(\cL1)$
and $(\cL2)$ of \cite[Exp.V, D\'ef.8.1.1]{SGA4-2}, hence theorem
\ref{th_gener-Beilinson} generalizes as well
\cite[Exp.V, Prop.8.1.12]{SGA4-2}.

(iii)\ \
In the situation of theorem \ref{th_gener-Beilinson}, one can
show that $(\cC',J')$ is not necessarily a $\sU$-site.
\end{remark}

We may now generalize lemma \ref{lem_cont-funct-site} to any
continuous functor between $\sU$-sites. Indeed, we notice
the following further application of theorem
\ref{th_gener-Beilinson}, which appears in
\cite[Exp.III, Th.4.1]{SGA4-1}.

\begin{proposition}\label{prop_comparis-lemma}
Let $C:=(\cC,J)$ be a $\sU$-site, $G$ a small topologically
generating family for $C$. Denote by $\cG$ the full subcategory
of\/ $\cC$ with $\Ob(\cG)=G$, and endow $\cG$ with the topology
$J_\cG$ induced by $J$ via the inclusion functor $u:\cG\to\cC$. Then :
\begin{enumerate}
\item
$u$ is cocontinuous for the topologies $J,J_\cG$ and is
a morphism of sites $C\to(\cG,J_\cJ)$.
\item
The induced functor $\tilde u_*:C^\sim\to(\cG,J_\cG)^\sim$
is an equivalence.
\item
A family $(\psi_i:Y_i\to Y~|~i\in I)$ of morphisms of\/ $\cG$
generates a covering sieve of $J_\cG$ if and only if it generates
a covering sieve of $J$.
\end{enumerate}
\end{proposition}
\begin{proof} We apply the criterion of theorem
\ref{th_gener-Beilinson}, and taking into account remark
\ref{rem_gener-Beilinson} we are reduced to checking condition
(a) of the theorem. The latter holds by definition of
topologically generating family.
\end{proof}

\begin{corollary}\label{cor_two-U-sites}
Let $C:=(\cC,J)$ and $C':=(\cC',J')$ be two $\sU$-sites,
$g:\cC\to\cC'$ a functor, and $\sV,\sV'$ a pair of universes
with $\sU\subset\sV\subset\sV'$. We have :
\begin{enumerate}
\item
If $g$ is continuous, the following holds :
\begin{enumerate}
\item
The functor $\tilde g_{\sV*}:C^{\prime\sim}_\sV\to C_\sV^\sim$
admits a left adjoint
$\tilde g{}_\sV^*:C_\sV^\sim\to C_\sV^{\prime\sim}$.
\item
There are essentially commutative diagrams of categories :
$$
\xymatrix{
C_\sV^\sim \ar[r]^-{\tilde g{}_\sV^*} \ar[d]_i &
C_\sV^{\prime\sim} \ar[d]^{i'} &
\cC \ar[r]^-g \ar[d]_{h^a_\cC} & \cC' \ar[d]^{h^a_{\cC'}} \\
C_{\sV'}^\sim \ar[r]^-{\tilde g{}_{\sV'}^*} &
C_{\sV'}^{\prime\sim} & C^\sim_\sU \ar[r]^-{\tilde g{}_\sU^*}
 & C_\sU^{\prime\sim}
}$$
where $i$ and $i'$ are the inclusion functors.
\item
Suppose moreover that\/ $\cC$ is finitely complete and
$g$ is left exact. Then $\tilde g{}_\sV^*$ is exact.
\end{enumerate}
\item
If $g$ is cocontinuous, the following holds :
\begin{enumerate}
\item
The functor $\breve g{}_\sV^*:C_\sV^{\prime\sim}\to C_\sV^\sim$ admits
a right adjoint $\breve g_{\sV*}:C^\sim_\sV\to C_\sV^{\prime\sim}$.
\item
There is an essentially commutative diagram of categories :
$$
\xymatrix{
C_\sV^\sim \ar[r]^-{\breve g_{\sV*}} \ar[d] &
C_\sV^{\prime\sim} \ar[d] \\
C_{\sV'}^\sim \ar[r]^-{\breve g_{\sV'*}} &
C_{\sV'}^{\prime\sim}
}$$
whose vertical arrows are the inclusion functors.
\end{enumerate}
\item
If $g$ is continuous and cocontinuous, we have
natural isomorphisms of functors :
$$
\tilde g_{\sV*}\isom\breve g{}^*_\sV
\ :\ C'^\sim_\sV\to C^\sim_\sV
\qquad\qquad
(-)^a\circ g^\wedge_\sV\isom\tilde g_{\sV*}\circ(-)^a
\ :\ \cC'^\wedge_\sV\to C^\sim_\sV.
$$
\end{enumerate}
\end{corollary}
\begin{proof} (i.a): We choose a small topologically generating
family $G$ for $C$, and define the site $(\cG,J_\cG)$ and the
continuous functor $u:(\cG,J_\cG)\to C$ as in proposition
\ref{prop_comparis-lemma}. By applying lemma
\ref{lem_cont-funct-site}(i) to the continuous functor
$v:=g\circ u$, we deduce that
$\tilde v_{\sV*}=\tilde u_{\sV*}\circ\tilde g_{\sV*}$
admits a left adjoint. Then the assertion follows from
proposition \ref{prop_comparis-lemma}(ii). 

(i.b): More precisely, 
$\tilde g{}^*_\sV=\tilde v{}^*_\sV\circ\tilde u_{\sV*}$. Thus,
the essential commutativity of the left diagram follows from
\eqref{subsec_adj-universes}. Likewise, the essential
commutativity of the right diagram follows from that of
the left diagram (with $\sV:=\sU$ and $\sV'$ large enough
so that $\cC$ is $\sV'$-small) together with lemma
\ref{lem_cont-funct-site}(ii).

(i.c) holds by example \ref{ex_simple-case}.

(ii.a): Let $u$ and $h$ be as in the foregoing; from
proposition \ref{prop_comparis-lemma}(i) we deduce that
$h$ is cocontinuous, hence
$\breve h{}^*_\sV=\breve u{}^*_\sV\circ\breve g{}^*_\sV$
admits a right adjoint; however
$\breve u{}^*_\sV=\tilde u_{\sV*}$ (lemma \ref{lem_improve}(i)),
and the latter is an equivalence (proposition
\ref{prop_comparis-lemma}(ii)), whence the contention.

(ii.b): More precisely,
$\breve g_{\sV*}=\breve h_{\sV*}\circ\tilde u{}^*_\sV$,
hence the assertion follows from \eqref{subsec_adj-universes}.

(iii): In view of (i.b) and (ii.b), in order to prove the
assertion, we may assume that $\cC$ and $\cC'$ are $\sV$-small,
and this case is already covered by lemma \ref{lem_improve}(i,ii).
\end{proof}

The last result of this section will show the representability
of the $2$-limit of any (small) cofiltered system of lex-sites.
The proof shall use the following :

\begin{lemma}\label{lem_filt-colims-lex-cats}
Let $I$ be a small filtered category, $\sU'$ a universe, and
$$
\cC_\bullet:I\to\sU'\tdu\bCat
\qquad
i\mapsto\cC_i
\qquad
(\phi:i\to j)\mapsto(\cC_\phi:\cC_i\to\cC_j)
$$
any pseudo-functor such that $\cC_i$ is finitely complete for
every $i\in\Ob(I)$ and $\cC_\phi$ is a left exact functor
for every morphism $\phi$ of $I$. Let also $(\cC,\pi_\bullet)$
be a $2$-colimit of\/ $\cC_\bullet$. Then $\cC$ is finitely complete
and $\pi_i:\cC_i\to\cC$ is a left exact functor for every
$i\in\Ob(I)$.
\end{lemma}
\begin{proof} Notice that the assertions depend only on
the equivalence class of the category $\cC$; hence, we
may suppose that $(\cC,\pi_\bullet)$ is the strong $2$-colimit
of $\cC_\bullet$ described explicitly by example
\ref{ex_filter-2-colim-in-Cat}(i), {\em i.e.}
$\cC=\cFib(\cC_\bullet)[\Sigma^{-1}]$, where $\Sigma$ is
the set of cartesian morphisms of $\cFib(\cC_\bullet)$;
then $\pi_i$ is induced by the inclusion functor of the
fibre category $F^{-1}(i)=\cC_i$ into $\cFib(\cC_\bullet)$,
where $F:\cFib(\cC_\bullet)\to I^o$ is the natural projection.
By proposition \ref{prop_complete-criteria}(i)), in order
to check that $\cC$ is finitely complete it suffices to
show the following two claims :

\begin{claim}\label{cl_fin-prods-ok}
All finite products are representable in $\cC_\bullet$.
\end{claim}
\begin{pfclaim} By a simple induction we are easily reduced to
showing that the product of two objects $(i^o_1,X_1),(i^o_2,X_2)$
is representable (here we have $i^o_1,i^o_2\in\Ob(I^o)$ and
$X_t\in\Ob(\cC_{i_t})$ for $t=1,2$ : see \eqref{subsec_fib-from-pseudo}).
Since $I^o$ is cofiltered, we may find $i\in\Ob(I)$ and
morphisms $\phi^o_t:i^o\to i^o_t$ for $t=1,2$; then $(i^o_t,X_t)$
is isomorphic to $(i,\cC_{\phi_t}X_t)$ in $\cC$ for $t=1,2$, so
we may assume that $i=i_1=i_2$. Since $\cC_i$ is finitely
complete, the product $X_1\times X_2$ is representable in
$\cC_i$, say by an object $P$, and let $(p_t:P\to X_t~|~t=1,2)$
be a universal cone. Let us check that
$$
(q_t:=(\one_{i^o},p_t):(i^o,P)\to(i^o,X_t)~|~t=1,2)
$$
is a universal cone in $\cC$. Indeed, let $(j^o,Y)$ be any
object of $\cC$ and $(r_t:(j^o,Y)\to(i^o,X_t)~|~t=1,2)$ a pair
of morphisms of $\cC$; we need to show that there exists a
unique morphism
\set\begin{equation}\label{eq_unique-r}
r:(j^o,Y)\to(i^o,P)
\qquad\text{in $\cC$ such that}\qquad
q_t\circ r=r_t
\qquad\text{for $t=1,2$}.
\end{equation}
By example \ref{ex_filter-2-colim-in-Cat}(ii) we can write
$r_t=g_t\circ s^{-1}_t$, where $s_t:(j^o_t,Y_t)\to(j^o,Y)$
and $g_t:(j^o_t,Y_t)\to(i^o,X_t)$ are morphisms of
$\cFib(\cC_\bullet)$, and $s_t$ is cartesian, for $t=1,2$.
By (CF3) of definition \ref{def_right-calculus}(i) we may
then find an object $(j'^o,Y')$ of $\cC$ and morphisms
$s'_t:(j'^o,Y')\to(j^o_t,Y_t)$ of $\cFib(\cC_\bullet)$ such that
$s'':=s_1\circ s'_1=s_2\circ s'_2$ and such that $s'_1$ is
cartesian; therefore $s''$ is an isomorphism in $\cC$, so
it suffices to check that there exists a unique morphism
$r':(j'^o,Y')\to(i^o,P)$ in $\cC$ with $q_t\circ r'=r_t\circ s''$
for $t=1,2$. Thus, we may assume that $r_t$ is the class of a
morphism $(\rho^o_t,f_t):(j^o,Y)\to(i^o,X_t)$ for $t=1,2$
(notation of \eqref{subsec_fib-from-pseudo}, so here
$\rho_t:i\to j$ is a morphism of $I$ and $f_t:Y\to\cC_{\rho_t}X_t$
is a morphism of $\cC_j$ for $t=1,2$). Since $I$ is filtered, we
may then find a morphism $\rho':j\to j'$ in $I$ such that
$\rho'\circ\rho_1=\rho'\circ\rho_2$, and since
$(\rho'^o,\one_{\cC_\rho}Y):(j'^o,\cC_\rho Y)\to(j^o,Y)$ is an
isomorphism in $\cC$, we may further replace $(\rho^o_t,f_t)$
by its composition with $(\rho'^o,\one_{\cC_\rho}Y)$ for $t=1,2$,
and assume as well that $\rho:=\rho_1=\rho_2$. In this situation
set $P'=\cC_\rho P$ and $X'_t:=\cC_\rho X_t$, $p'_t:=\cC_\rho(p_t)$
for $t=1,2$; we notice that $r_t$ factors through a morphism
$r'_t:=(\one_{j^o},f'_t):(j^o,Y)\to(j^o,X'_t)$ and the cartesian
morphism $(\rho,\one_{X'_t}):(j^o,X'_t)\to(i^o,X_t)$, and moreover
$$
(\rho,\one_{X'_t})\circ(\one_{j^o},p'_t)=
(\one_{i^o},p_t)\circ(\rho,\one_{P'})
\qquad
\text{in $\cFib(\cC_\bullet)$ for $t=1,2$}.
$$
Since $\cC_\rho$ is left exact, the cone $(p'_t:P'\to X'_t~|~t=1,2)$
is universal in $\cC_j$, so there exists a unique morphism
$g':Y\to P'$ in $\cC_j$ such that $p'_t\circ g'=f'_t$ for $t=1,2$.
It follows that $r:=(\rho,g'):(j^o,Y)\to(i^o,P)$ fulfills the
condition of \eqref{eq_unique-r}. It remains to check the
uniqueness of $r$. Thus, let $r':(j^o,Y)\to(i^o,P)$ be another
morphism of $\cC$ that verifies the identities of
\eqref{eq_unique-r}. We have $r'=g'\circ s'^{-1}$ for some
morphisms $g'$ and $s'$ of $\cFib(\cC_\bullet)$ with $s'$
cartesian, and we are reduced to checking that $g'=r\circ s'$.
We may then assume that $r'$ is the class of a morphism
$(\rho'^o,g'')$ of $\cFib(\cC_\bullet)$. Moreover, the
identities \eqref{eq_unique-r} for $r'$ mean that there
exist cartesian morphisms $v_t:(k^o_t,Y'_t)\to(j^o,Y)$ such
that $q_t\circ r'\circ v_t=r_t\circ v_t$ in $\cFib(\cC_\bullet)$
for $t=1,2$. By (CF3) we may find an object $(k^o,Y'')$ of
$\cFib(\cC_\bullet)$ and morphisms $v'_t:(k^o,Y'')\to(k^o_t,Y'_t)$
in $\cFib(\cC_\bullet)$ with $v'_1$ cartesian, and such that
$v'':=v_1\circ v'_1=v_2\circ v'_2$. Since $v''$ is an isomorphism
in $\cC$, it suffices to check that $r\circ v''=r'\circ v''$,
and we may therefore assume that the identities \eqref{eq_unique-r}
hold already in the category $\cFib(\cC_\bullet)$, for both $r$
and $r'$.

We may next find a morphism $\nu:j'\to j$ in $I$ such that
$\nu\circ\rho=\nu\circ\rho'$, and it suffices to check that
$(\nu^o,\one_{\cC_\nu Y})\circ r=(\nu^o,\one_{\cC_\nu Y})\circ r'$.
We may therefore assume that $\rho=\rho'$. Then both $r$ and
$r'$ factor uniquely through morphisms
$(\one_{j^o},u),(\one_{j^o},u'):(j^o,Y)\to(j^o,P')$ and the
cartesian morphism $(\rho^o,\one_{P'}):(j^o,P')\to(i^o,P)$,
and we have
$$
(\one_{j^o},p'_t)\circ(\one_{j^o},u)=r'_t=
(\one_{j^o},p'_t)\circ(\one_{j^o},u')
\qquad
\text{in $\cFib(\cC_\bullet)$ for $t=1,2$}.
$$
Since $P'$ represents $X'_1\times X'_2$ in $\cC_j$,
we conclude that $u=u'$, whence $r=r'$ as required.
\end{pfclaim}

\begin{claim}\label{cl_equals-ok}
All equalizers are representable in $\cC$.
\end{claim}
\begin{pfclaim} Let $(i^o_1,X_1),(i^o_2,X_2)$ be objects
of $\cC$ and $f_1,f_2:(i^o_1,X_1)\to(i^o_2,X_2)$ morphisms of
$\cC$. We need to represent the equalizer of $f_1$ and $f_2$, and
arguing as in the foregoing we may assume that $i:=i_1=i_2$,
and also that $f_t$ is the class of a morphism $(\phi_t,g_t)$
of $\cFib(\cC_\bullet)$ for $t=1,2$. Then we may further reduce
to the case where $\phi:=\phi_1=\phi_2$. After replacing
$(i^o_2,X_2)$ by $(i^o_1,\cC_\phi X_2)$, we may then assume
that $f_t=(\one_{i^o},g_t):(i^o,X_1)\to(i^o,X_2)$ for $t=1,2$.
By assumption, the equalizer of $g_1$ and $g_2$ in $\cC_i$
is representable by some $E\in\Ob(\cC_i)$, and the universal
cone for this equalizer amounts to a morphism $u:E\to X_1$
in $\cC_i$ such that $g_1\circ u=g_2\circ u$. We claim that
$(i^o,E)$ represents the equalizer of $f_1$ and $f_2$, and
that the morphism $(\one_{i^o},u):(i^o,E)\to(i^o,X_1)$ yields
a universal cone for this equalizer. Indeed, let $(j^o,Y)$
be an object of $\cC$ and $h:(j^o,Y)\to(i^o,X_1)$ a morphism
in $\cC$ such that $f_1\circ h=f_2\circ h$; we need to show
that $h$ factors uniquely through $u$ and a morphism
$h':(j^o,Y)\to(i^o,E)$. We reduce easily to the case where
$h$ is the class of a morphism $(\psi^o,k):(j^o,Y)\to(i^o,X_1)$
of $\cFib(\cC_\bullet)$ such that
$(\one_{i^o},g_1)\circ(\psi^o,k)=(\one_{i^o},g_2)\circ(\psi^o,k)$
in $\cFib(\cC_\bullet)$. In this situation, set $E':=\cC_\psi E$,
$u':=\cC_\psi u$, $X'_t:=\cC_\psi X_t$ and $g'_t:=\cC_\psi g_t$ for
$t=1,2$. Then $(\psi^o,k)$ factors uniquely through the morphism
$(\one_{j^o},k):(j^o,Y)\to(j^o,X'_1)$ and the cartesian morphism
$(\psi^o,\one_{X'_1}):(j^o,X'_1)\to(i^o,X_1)$, and we have
$(\one_{j^o},g'_1)\circ(\one_{j^o},k)=
(\one_{j^o},g'_2)\circ(\one_{j^o},k)$. Since $\cC_\psi$ is
left exact, $E'$ represents the equalizer of $g'_1$ and $g'_2$
in $\cC_j$, so there exists a unique morphism $k':Y\to E'$ in
$\cC_j$ such that $u'\circ k'=k$. Then $h:=(\psi,u')$ is the
sought morphism. The verification of the uniqueness of $h$ is
analogous to that of the corresponding assertion in the proof
of claim \ref{cl_fin-prods-ok} : the details shall be left to
the reader.
\end{pfclaim}

Lastly, the constructions in the proof of claims
\ref{cl_fin-prods-ok} and \ref{cl_equals-ok} show that the
functors $\pi_i$ commute with finite products and equalizers,
so they are left exact, by proposition \ref{prop_commute-criteria}(i).
\end{proof}

\sset\subsubsection{}\label{subsec_cofilt-lims-lexsites}
Let now $I$ be a small cofiltered category and
$C_\bullet:I\to(\sU,\sU')\tdu\lex.\Site$ any pseudo-functor.
For every $i\in\Ob(I)$, say that $C_i=(\cC_i,J_i)$, and denote by
$\cC_\bullet:I^o\to\sU'\tdu\bCat$ the composition of $C^o_\bullet$
with the forgetful functor
${}^o(\sU,\sU')\tdu\lex.\Site^o\to\sU'\tdu\bCat$. We let
$(\cC,\pi_\bullet)$ be a strong $2$-colimit of the pseudo-functor
$\cC_\bullet$ (theorem \ref{th_bCat-cpt}), and we endow $\cC$ with
the coarsest topology $J$ such that all the functors
$\pi_i:\cC_i\to\cC$ are continuous for the topologies $J_i$ and
$J$. In light of example \ref{ex_simple-case} and lemma
\ref{lem_filt-colims-lex-cats}, we see that $C:=(\cC,J)$ is a
lex-site, and we obtain a well defined pseudo-cone
$\pi_\bullet^o:\sF_C\Rightarrow C_\bullet$.

\begin{proposition}\label{prop_2-lims-lex-sites}
In the situation of \eqref{subsec_cofilt-lims-lexsites},
let $F\in\Ob(\cC^\wedge_\sU)$ be any presheaf. We have :
\begin{enumerate}
\item
The pair $(C,\pi^o_\bullet)$ is a strong $2$-limit of\/ $C_\bullet$.
\item
$F$ is a sheaf on $C$ if and only if $\pi^\wedge_i(F)$ is a sheaf
on $C_i$ for every $i\in\Ob(I)$.
\end{enumerate}
\end{proposition}
\begin{proof}(i): Let $D:=(\cD,J_\cD)$ be any finitely complete
$\sU$-site, and $\phi_\bullet:\sF_D\Rightarrow C_\bullet$ any
pseudo-cone. Hence for each $i\in\Ob(I)$ the functor
$\phi_i:\cC_i\to\cD$ is left exact and continuous for the
topologies $J_i$ and $J_\cD$. There results a pseudo-cocone
$\phi^a_\bullet:\cC_\bullet\Rightarrow\sF_\cD$, whence a unique
functor $\phi:\cC\to\cD$ such that $\sF_\phi\odot\pi^o_\bullet=\phi^a$.
The assertion then follows immediately from :

\begin{claim}(i)\ \ The functor $\phi$ is left exact.

(ii)\ \ The functor $\phi$ is a morphism of sites $C\to D$.

(iii)\ \ $C$ is a $\sU$-site.
\end{claim}
\begin{pfclaim}(i): It suffices to check that $\phi$
commutes with all finite products and all equalizers of $\cC$
(proposition \ref{prop_commute-criteria}(i)). Thus, let
$X_1,X_2$ be two objects of $\cC$, and $P$ another object
that represents $X_1\times X_2$; pick also a universal
cone for this product, given by a pair of morphisms
$p_\bullet:=(p_t:P\to X_t~|~t=1,2)$. By inspecting the proof
of claim \ref{cl_fin-prods-ok}, we see that there exist
$i\in\Ob(I)$ and objects $P',X'_1,X'_2\in\Ob(\cC_i)$ such
that $P'$ represents $X'_1\times X'_2$ and the cone $p_\bullet$
is isomorphic to $\cC_i*p'_\bullet$, where
$p'_\bullet:=(p'_t:P'\to X'_t~|~t=1,2)$ is a universal cone in
$\cC_i$. Then $\phi*p_\bullet$ is isomorphic to $\phi_i*p'_\bullet$,
and the latter is universal, since $\phi_i$ is left exact.
This shows that $\phi$ commutes with binary products, and hence
with all finite products. Similarly, by inspecting the proof
of claim \ref{cl_equals-ok} we see that the universal cone
of every equalizer in $\cC$ is isomorphic to $\cC_i*\tau$,
for some $i\in\Ob(I)$ and with $\tau$ the universal cone of
some equalizer in $\cC_i$; then again we deduce easily that
$\phi$ commutes with equalizers.

(ii): In view of (i) and example \ref{ex_simple-case}, it
suffices to show that $\phi$ is continuous for the topologies
$J$ and $J_\cD$. To this aim, let $J'$ be the topology on $\cC$
induced by $J_\cD$ via $\phi$; it suffices to check that
$J\subset J'$. However, from lemma \ref{lem_crit-continuity}
we see that $J$ is the coarsest topology on $\cC$ such that
for every $i\in\Ob(I)$, every $X\in\Ob(\cC_i)$ and every
covering family $(g_\lambda:X_\lambda\to X~|~\lambda\in\Lambda)$
for the topology $J_i$, the sieve generated by the family
$(\pi_i(g_\lambda)~|~\lambda\in\Lambda)$ lies in $J(\pi_i(X))$.
Thus, consider any such covering family
$(g_\lambda~|~\lambda\in\Lambda)$; according to lemma
\ref{lem_induced-top}(ii.b) we are further reduced to checking
that the family $(\phi\circ\pi_i(g_\lambda)~|~\lambda\in\Lambda)$
covers $\phi\circ\pi_i(X)=\phi_i(X)$ for the topology $\cD$.
Since $\phi\circ\pi_i$ is continuous and $\cC_i$ is finitely
complete, the latter follows again from lemma
\ref{lem_crit-continuity}.

(iii): For every $i\in\Ob(I)$ pick a small topologically
generating family $G_i\subset\Ob(\cC_i)$ for the site $C_i$.
It suffices to show that $G:=\bigcup_{i\in\Ob(I)}\pi_i(G_i)$
is a topologically generating family for the site $C$.
Thus, let $X\in\Ob(\cC)$ be any object; by inspecting the
construction of $\cC$ we see that there exist $i\in\Ob(I)$
and $X_i\in\Ob(\cC_i)$ such that $X=\pi_i(X_i)$, and by
assumption the family $\cF:=\bigcup_{Y\in G_i}\Hom_{\cC_i}(Y,X_i)$
covers $X_i$ for the topology $J_i$. It follows that the family
$(\pi_i(g)~|~g\in\cF)\subset\bigcup_{Z\in G}\Hom_\cC(Z,X)$ covers
$X$ for the topology $J$, whence the contention.
\end{pfclaim}

(ii): The condition is obviously necessary. For the converse,
we consider the finest topology $J''$ on $\cC$ such that $F$
is a sheaf on the site $(\cC,J'')$ (see remark
\ref{rem_sheaves}(iii)); it then suffices to show that
$J\subset J''$. In view of remark \ref{rem_sheaves}(iii), we
come down to checking the following. For every $i\in\Ob(I)$,
every $X\in\Ob(\cC_i)$, every covering family
$(g_\lambda:X_\lambda\to X~|~\lambda\in\Lambda)$
for the topology $J_i$ and every morphism $f:Y\to\pi_iX$
in $\cC$, the family
$$
(Y\times_{\pi_iX}\pi_i(g_\lambda):
Y\times_{\pi_iX}\pi_iX_\lambda\to Y~|~\lambda\in\Lambda)
$$
generates a sieve of $2$-descent for the fibration
$\cFib(F)\to\cC$ of \eqref{subsec_fibred-cats-II}.
But recall that $\cC=\cFib(\cC_\bullet)[\Sigma^{-1}]$, where
$\Sigma$ denotes the set of cartesian morphisms of
$\cFib(\cC_\bullet)$. With this notation, we have $\pi_iX=(i^o,X)$
and $Y=(j^o,Y')=\pi_jY'$ for some $j\in\Ob(I)$ and $Y'\in\Ob(\cC_j)$.
We are then easily reduced to the case where $f$ is the class
of a morphism $(\phi^o,t):(j^o,Y')\to(i^o,X)$ of $\cFib(\cC_\bullet)$,
for some morphism $\phi:i\to j$ in $I$ and $t:Y'\to\cC_\phi X$
in $\cC_j$. Notice now the commutative diagrams :
$$
{\diagram
\pi_j(\cC_\phi X_\lambda) \ar[rr]^-{\pi_j(\cC_\phi g_\lambda)}
\ar[d]_{(\phi,\one_{\cC_\phi X_\lambda})} & &
\pi_j(\cC_\phi X) \ar[d]^{(\phi,\one_{\cC_\phi X})} \\
\pi_i(X_\lambda) \ar[rr]^-{\pi_i(g_\lambda)} & & \pi_i(X)
\enddiagram}
\qquad
\text{for every $\lambda\in\Lambda$}
$$
whose vertical arrows are cartesian morphisms of
$\cFib(\cC_\bullet)$, so they are isomorphisms in $\cC$.
Moreover, the continuity of $\cC_\phi$ implies that the family
$(\cC_\phi(g_\lambda)~|~\lambda\in\Lambda)$ covers $\cC_\phi X$
for the topology $J_j$. Furthermore, we have
$(\phi^o,t)=(\one_{j^o},t)\circ(\phi,\one_{\cC_\phi X})$. Thus,
we may replace $(g_\lambda~|~\lambda\in\Lambda)$ by
$(\cC_\phi(g_\lambda)~|~\lambda\in\Lambda)$ and assume that
$i=j$ and $\phi=\one_i$. In this situation, since $\pi_i$
is left exact, we have natural identifications for every
$\lambda,\mu\in\Lambda$ :
$$
Y_\lambda:=Y\times_{\pi_iX}\pi_iX_\lambda\isom\pi_i(Y'\times_XX_\lambda)
\qquad
Y_\lambda\times_YY_\mu\isom\pi_i(Y'\times_XX_\lambda\times_XX_\mu).
$$
So finally we are reduced to checking that the natural map
from $F(\pi_iY')$ to the equalizers of the restriction maps
$$
\prod_{\lambda\in\Lambda}\xymatrix{ F\circ\pi_i(Y'\times_XX_\lambda)
\ar@<.5ex>[r] \ar@<-.5ex>[r] &}
\prod_{\lambda,\mu\in\Lambda}F\circ\pi_i(Y'\times_XX_\lambda\times_XX_\mu)
$$
is a bijection. But this is clear, since by assumption $F\circ\pi_i$
is a sheaf on $C_i$.
\end{proof}

\begin{remark}\label{rem_2-lims-of-sites}
In the situation of \eqref{subsec_cofilt-lims-lexsites}, we also
claim that $C^\sim$ represents the $2$-limit of the induced
pseudo-functor
$$
C^\sim_\bullet:I\to\sU'\tdu\bCat
\qquad
i\mapsto C_{i,\sU}^\sim
\qquad
(\phi:i\to j)\mapsto(C_\phi^\sim:C^\sim_{i,\sU}\to C_{j,\sU}^\sim).
$$
Indeed, on the one hand, the natural functor (see definition
\ref{def_pseudo-lim}(ii))
$$
\sPsFun(\cC,\Set^o)\to\sPsNat(\cC_\bullet,\sF_{\Set^o})
\qquad
(G:\cC\to\Set^o)\mapsto\sF_G\odot\pi_\bullet
$$
is an equivalence; but we have $\sPsFun(\cC,\Set^o)=\cC^\wedge$,
and a pseudo-cocone $G_\bullet:\cC_\bullet\Rightarrow\sF_{\Set^o}$
is the datum of a system of presheaves $G_i$ on $\cC_i$ for
every $i\in\Ob(I)$, with isomorphisms
$$
\tau^G_\phi:G_i\isom\cC_\phi^\wedge G_j
\qquad
\text{in $\cC^\wedge_i$ for every morphism $\phi:i\to j$ in $I$}
$$
fulfilling the usual coherence axioms.  Under these
identifications, a presheaf $G$ on $\cC$ corresponds to the
datum $(G_\bullet,\tau^G_\bullet)$ with $G_i:=\pi^\wedge_iG$ for
every $i\in\Ob(I)$, and with
$$
\tau^G_\phi:=(\tau^\pi_\phi)^\wedge_G:\pi^\wedge_iG\isom
\cC^\wedge_\phi\circ\pi^\wedge_jG
$$
where
$(\tau^\pi_\phi)^\wedge:\pi^\wedge_i\isom\cC^\wedge_\phi\circ\pi^\wedge_j$
is the isomorphism of functors induced by the coherence constraint
$\tau^\pi_\phi:\pi_j\isom\cC_\phi\circ\pi_i$ of the
universal pseudo-cone $\pi_\bullet:\sF_\cC\Rightarrow\cC_\bullet$,
for every morphism $\phi:i\to j$ in $I$. On the other hand, by
proposition \ref{prop_2-lims-lex-sites}(ii), the presheaf $G$ is
a sheaf on $C$ if and only if $G_i$ is a presheaf on $C_i$ for
every $i\in\Ob(I)$. Summing up, we obtain a natural equivalence
between $C^\sim$ and the category of pairs $(G_\bullet,\tau^G_\bullet)$
as above, such that $G_i$ is a sheaf on $C_i$ for every $i\in\Ob(I)$.
But the latter is just a description of the category
$$
\sCart_{I^o}(I^o,\cFib(C_\bullet^\sim))
$$
which represents the $2$-limit of $C_\bullet^\sim$, by the proof
of theorem \ref{th_bCat-cpt}.
\end{remark}

\subsection{Topoi}
Let $\sU,\sU'$ be a pair of universes with $\sU\in\sU'$.
A {\em $(\sU,\sU')$-topos\/} is a $\sU'$-small category
with $\sU$-small $\Hom$-sets which is equivalent to the category
of sheaves on a $\sU$-small site. We shall often omit the explicit
mention of $\sU'$, and call ``$\sU$-topos'' such a category, or
even just ``topos'' unless this may give rise to ambiguities.

\begin{remark}\label{rem_summarized}
Let $T$ be any $\sU$-topos.

(i)\ \
Remark \ref{rem_rep-and-sheafify}(i,ii) can be summarized by
saying that $T$ is a complete and cocomplete, well-powered
and co-well-powered category. Morever, every epimorphism in
$T$ is universal effective, every colimit is universal, and
all filtered colimits in $T$ commute with finite limits.

(ii)\ \
It follows from (i) that a small family $(F_i\to F~|~i\in I)$ of
morphisms of $T$ is a covering family for the canonical topology
if and only if the induced morphism $\amalg_{i\in I}F_i\to F$ is
an epimorphism in $T$.

(iii)\ \
Say that $T=C^\sim$, for a small site $C:=(\cC,J)$; then the set
$\{h^a_X~|~X\in\Ob(\cC)\}$ is a small topologically generating
family for the canonical site $\Can(T)$ (see remark
\ref{rem_sheaves}(v)). Indeed, by virtue of remark
\ref{rem_rep-and-sheafify}(iii) we may find for every
$F\in\Ob(T)$ a family $(X_i~|~i\in I)$ of objects of $\cC$ with
an epimorphism $\amalg_{i\in I}h^a_{X_i}\to F$ in the category
$\cC^\sim$, so the resulting family $(h^a_{X_i}\to F~|~i\in I)$
covers $F$ for the canonical topology, by (ii). Especially,
$\Can(T)$ is a $\sU$-site.

(iv)\ \
On the other hand, by proposition \ref{prop_comparis-lemma}(ii),
if $C$ is a $\sU$-site, then $C^\sim$ is isomorphic to a $\sU$-topos.
\end{remark}

\begin{theorem}\label{th_canon-topos}
Let $C:=(\cC,J)$ be a site such that $\cC$ has small
$\Hom$-sets. We have :
\begin{enumerate}
\item
The Yoneda embedding $\cC\to\cC^\wedge$ is a morphism of sites
$h_\cC:C^\wedge\to C$ (notation of remark
{\em\ref{rem_topol-on-presheaves}(ii)}) and is cocontinuous for
the topologies $J$ and $J^\wedge$. Moreover the topology $J$ is
induced by $J^\wedge$ via the functor $h_\cC$.
\item
Suppose that $C$ is a $\sU$-site. Then also the forgetful
functor $i:C^\sim\to\cC^\wedge$, its left adjoint
$(-)^a:\cC^\wedge\to C^\sim$ and the Yoneda embedding
$h^a_C:\cC\to C^\sim$ are morphisms of sites
$C^\wedge\to\Can(C^\sim)$, $\Can(C^\sim)\to C^\wedge$
and $\Can(C^\sim)\to C$ respectively.
\item
Moreover, in the situation of {\em(ii)} the functors $(-)^a$ and
$h^a_C$ are cocontinuous for the topologies $J$, $J^\wedge$
and $\Can_{C^\sim}$, and the topologies $J$ and $J^\wedge$ are
induced by the canonical topology of $C^\sim$, via the functors
$h^a_C:\cC\to C^\sim$ and $(-)^a$ respectively. Furthermore,
for every universe $\sV$ containing $\sU$, the functors $i$,
$h_\cC$, $h^a_C$, $(-)^a$ induce equivalences :
$$
\xymatrix{ \Can(C^\sim)^\sim_\sV \ar@<.5ex>[rr]^-{((-)^a)^\sim_*}
\ar[rd]_{\tilde h{}^a_{C*}} & &
(C^\wedge)^\sim_\sV \ar@<.5ex>[ll]^-{\tilde\imath_*} \ar[ld]^{\tilde h_{\cC*}} \\
& C^\sim_\sV.
}$$
\item
For every topos $T$, the Yoneda embedding induces an equivalence
$$
h_T:T\to\Can(T)^\sim
\qquad
F\mapsto h_F.
$$
\end{enumerate}
\end{theorem}
\begin{proof}(i): Pick a universe $\sV$ containing $\sU$, and
such that $C^\wedge$ is a $\sV$-site; denote by $J'$ the topology
on $\cC$ induced by $J^\wedge$ on $\cC$ via $h_\cC$. Notice that
condition (a) of theorem \ref{th_gener-Beilinson} holds for the
functor $h_\cC$, by virtue of remark \ref{rem_summarized}(iii),
and conditions (b) and (c) hold as well, since $h_\cC$ is fully
faithful (remark \ref{rem_gener-Beilinson}(i)). By theorem
\ref{th_gener-Beilinson}(iii) it follows that a family of
morphisms $(Y_i\to Y~|~i\in I)$ in $\cC$ covers $Y$ in the topology
$J'$ if and only if the family $(h_{Y_i}\to h_Y~|~i\in I)$ covers
$h_Y$ in the topology $J^\wedge$. Combining with remark
\ref{rem_iprippi}(i) we deduce that $J'=J$. Moreover, by
theorem \ref{th_gener-Beilinson}(ii), the functor
$(h_\cC)^\sim_{\sV*}:(C^\wedge)^\sim_\sV\to C^\sim_\sV$ is an equivalence,
so the same holds for its left adjoint $(h^*_\cC)^\sim_\sV$, and
especially, $h_\cC$ is a morphism of sites, as stated.

(ii): Let us remark :

\begin{claim}\label{cl_abc-again}
The functor $(-)^a:\cC^\wedge\to C^\sim$ satisfies conditions
(a),(b),(c) of theorem \ref{th_gener-Beilinson} for the
topology $\Can_{C^\sim}$.
\end{claim}
\begin{pfclaim} Indeed, to check condition (a)
it suffices to remark that $(iF)^a$ is isomorphic to $F$,
for every $F\in\Ob(C^\sim)$. Next, let $F,G\in\Ob(\cC^\wedge)$
and $\phi:F^a\to G^a$ a morphism in $C^\sim$; we may find a
family $(f_j:h_{X_j}\to F~|~j\in I)$ of morphisms in $\cC^\wedge$
such that the induced morphism $\mu:\amalg_{j\in I}h_{X_j}\to F$
is an epimorphism (lemma \ref{lem_lable}). The composition
$$
\phi_j:h_{X_j}\xrightarrow{f_j} F\to F^a\xrightarrow{\phi}G^a
\qquad
\text{for every $j\in J$}
$$
corresponds to a section $\sigma_j\in G^a(X_j)$. Then we may
find for every $j\in I$ a covering family
$(f_{jk}:X_{jk}\to X_j~|~k\in I_j)$ for the topology $J$,
and a section $\sigma_{jk}:G(X_{jk})$ for every $k\in I_j$,
whose image in $G^a(X_{jk})$ agrees with $G^a(f_{jk})(\sigma_j)$.
The section $\sigma_{jk}$ corresponds to a morphism
$\psi_{jk}:h_{X_{jk}}\to G$ in $\cC^\wedge$, whose composition
with the natural morphism $G\to G^a$ agrees with
$\phi_j\circ h_{f_{jk}}$. Moreover, for every $j\in J$ the
resulting morphism $\nu_j:\amalg_{k\in I_j}h_{X_{jk}}\to h_{X_j}$
of $\cC^\wedge$ induces an epimorphism $\nu^a$ in $C^\sim$
(remark \ref{rem_iprippi}(i) and corollary \ref{cor_bicover}(i));
likewise $\mu^a$ is an epimorphism in $C^\sim$, by virtue of
proposition \ref{prop_was-get-maddd}(v). Thus, the family
$(f_j^a\circ h^a_{f_{jk}}:h^a_{X_{jk}}\to F^a~|~j\in I,\ k\in I_j)$
covers $F^a$ for the topology $\Can_{C^\sim}$ (remark
\ref{rem_summarized}(ii)), and by construction we have
$\phi\circ(f_j\circ h_{f_{jk}})^a=\psi_{jk}^a$ for every $j\in I$
and every $k\in I_j$. This shows that condition (b) holds.
Lastly, let $\phi,\phi':F\to G$ be two morphisms of $\cC^\wedge$
such that $\phi^a=\phi'^a$, and choose again a family
$(f_j:h_{X_j}\to F~|~j\in I)$ as in the foregoing; each $f_j$
corresponds to a section $\sigma_j\in FX_j$, and by assumption
the images of $\phi_{X_j}(\sigma_j)$ and $\phi'_{X_j}(\sigma_j)$
agree in $G^a(X_j)$. Then we may find for every $j\in I$ a
covering family $(f_{jk}:X_{jk}\to X_j~|~k\in I_j)$ such that
$$
\phi_{X_{jk}}\circ(Ff_{jk})(\sigma_j)=(Gf_{jk})\circ\phi_{X_j}(\sigma_j)=
(Gf_{jk})\circ\phi'_{X_j}(\sigma_j)=\phi'_{X_{jk}}\circ(Ff_{jk})(\sigma_j)
$$
for every $k\in I_j$. Again, the family
$(f_j^a\circ h^a_{f_{jk}}:h^a_{X_{jk}}\to F^a~|~j\in I,\ k\in I_j)$
covers $F^a$ for the topology $\Can_{C^\sim}$, and
$\phi\circ f_j\circ h_{f_{jk}}=\phi'\circ f_j\circ h_{f_{jk}}$
for every $j\in I$ and every $k\in I_j$. This shows that also
condition (c) holds.
\end{pfclaim}

Let $J'$ be the topology on $\cC^\wedge$ induced by $\Can_{C^\sim}$
via the functor $(-)^a$; since $C^\sim$ is isomorphic to a $\sU$-site
(remark \ref{rem_summarized}(iv)), theorem \ref{th_gener-Beilinson}
and claim \ref{cl_abc-again} imply that
$(-)^a:(\cC^\wedge,J')\to\Can(C^\sim)$ is continuous and
cocontinuous, and a family $(F_j\to F~|~j\in I)$ of morphisms
in $\cC^\wedge$ generates a covering sieve for the topology $J'$
if and only if the family $(F^a_j\to F^a~|~j\in I)$ generates
a covering sieve for the topology $\Can_{C^\sim}$. But the latter
condition holds if and only if the family $(F_j\to F~|~j\in I)$
generates a covering sieve for the topology $J^\wedge$, by virtue
of remark \ref{rem_topol-on-presheaves}(i). In other words,
$J'=J^\wedge$, and thus theorem \ref{th_gener-Beilinson} also
says that $(-)^a$ induces an equivalence
$((-)^a)^\sim_{\sV*}:\Can(C^\sim)_\sV^\sim\isom(\cC^\wedge,J')_\sV^\sim$
for every universe $\sV$ containing $\sU$.
Then also the left adjoint of $((-)^a)^\sim_{\sV*}$ is an equivalence,
and especially $(-)^a:\Can(C^\sim)\to C^\wedge$ is a morphism of sites.
Lastly, since $(-)^a\circ i=\one_{C^\sim}$, we see that
$\tilde\imath_{\sV*}:(C^\wedge)_\sV^\sim\to\Can(C^\sim)_\sV^\sim$
is an equivalence as well, and so is its left adjoint (which is
$((-)^a)^\sim_{\sV*}$); especially, $i:C^\wedge\to\Can(C^\sim)$
is a morphism of sites.

\begin{claim}\label{cl_apply-th-to-Yoneda}
The Yoneda embedding $h^a_C:\cC\to C^\sim$ fulfills conditions
(a),(b),(c) of theorem \ref{th_gener-Beilinson}, for the
canonical topology on $C^\sim$.
\end{claim}
\begin{pfclaim} Indeed, condition (a) has already
been noticed in remark \ref{rem_summarized}(iii). Next, let
$X,Y\in\Ob(\cC)$ and $\phi:h^a_X\to h^a_Y$ any morphism of $C^\sim$;
then $\phi$ corresponds to a section $s\in h^a_Y(X)$, which
in turn is represented by the datum of a covering family
$(f_i:X_i\to X~|~i\in I)$ for the topology $J$ and a compatible
system of morphisms $(s_i:X_i\to Y~|~i\in I)$ in $\cC$. With
this notation, $h^a_{f_i}(s)$ is the image of $s$ in $h^a_Y(X_i)$,
so we get $h^a_{s_i}=\phi\circ h^a_{f_i}$ for every $i\in I$.
Lastly, the family $(h^a_{f_i}:h^a_{X_i}\to h^a_X~|~i\in I)$
covers $h^a_{X_i}$ in the canonical topology, by remarks
\ref{rem_summarized}(ii) and \ref{rem_iprippi}(i) and
corollary \ref{cor_bicover}(i).

To check condition (c), consider morphisms $\phi,\psi:X\to Y$
in $\cC$ such that $h^a_\phi=h^a_\psi$. The latter means that
the images $\bar\phi$ and $\bar\psi$ of $\phi$ and $\psi$
agree in $h^a_Y(X)$, so there exists a covering family
$(f_i:X_i\to X~|~i\in I)$ such that
$h^a_Y(f_i)(\bar\phi)=h^a_Y(f_i)(\bar\psi)$ for every $i\in I$.
This in turn translates as the identity
$\phi\circ f_i=\psi\circ f_i$ for every $i\in I$. To conclude,
we observe as in the foregoing that the family
$(h^a_{f_i}:h^a_{X_i}\to h^a_X~|~i\in I)$ covers $h^a_X$ in the
topology $\Can_{C^\sim}$.
\end{pfclaim}

Let $J'$ be the topology on $\cC$ induced by $\Can_{C^\sim}$ via
$h^a_\cC$; since we have already noticed that $C^\sim$ is isomorphic
to a $\sU$-site, claim \ref{cl_apply-th-to-Yoneda} and theorem
\ref{th_gener-Beilinson} imply that $h^a_\cC$ induces an equivalence
$H_\sV:=(h{}^a_C)^\sim_{\sV*}:\Can(C^\sim)_\sV^\sim\isom(\cC,J')^\sim_\sV$
for every universe $\sV$ containing $\sU$.
Moreover, arguing as in the foregoing  we deduce that $J'=J$.
Thus, we get an equivalence $H:\Can(C^\sim)^\sim\isom C^\sim$,
which as usual implies that $h^a_C:C^\sim\to C$ is a morphism of
sites, concluding the proof of (ii) and (iii). Lastly, by a simple
inspection we see that $H\circ h_{C^\sim}$ is naturally isomorphic
to $\one_{C^\sim}$, whence (iv).
\end{proof}

\begin{definition}\label{def_morph-of-topoi}
(i)\ \
Let $S$ and $T$ be two topoi. A {\em morphism of topoi\/}
$f:T\to S$ is a datum
$$
(f^*,f_*,\eta)
\qquad\text{where}\qquad
f_*:T\to S\qquad f^*:S\to T
$$
are two functors such that $f^*$ is left exact and left adjoint
to $f_*$, and $\eta:\one_S\Rightarrow f_*f^*$ is a unit of an
adjunction $\theta$ for the pair $(f^*,f_*)$ (these are sometimes
called {\em geometric morphisms\/} : see \cite[Def.2.12.1]{BorIII}).
We shall say that $\eta$ and $\theta$ are respectively
{\em the unit and the adjunction of $f$}.

(ii)\ \
Let $f:=(f^*,f_*,\eta_f):T''\to T'$ and $g:=(g^*,g_*,\eta_g):T'\to T$
be two morphisms of topoi. The composition $g\circ f$ is the
morphism
$$
(f^*\circ g^*,g_*\circ f_*,(g_**\eta_f*g^*)\odot\eta_g):T''\to T
$$
(see remark \ref{rem_adjoint-transf}(i)). To check the associativity
of this composition law, it suffices to notice that
$(g_**\eta_f*g^*)\odot\eta_g$ is the unit of the composition of the
adjunctions determined by $\eta_f$ and $\eta_g$, as defined in remark
\ref{rem_adjoint-transf}(i). For every topos $T$, let
$i_T:\one_T\isom\one_T$ be the identity automorphism of
$\one_T$; the datum
$$
(\one_T,\one_T,,i_T)
$$
yields a morphism of topoi $T\to T$ which is neutral for
left and right compositions relative to the foregoing
composition law.

(iii)\ \
Let $f,g:T\to S$ be two morphisms of topoi. A {\em natural
transformation\/} $\tau:f\Rightarrow g$ is just a natural
transformation of functors $\tau_*:f_*\Rightarrow g_*$.
Notice that, in view of remark \ref{rem_adjoint-transf}(ii), 
the datum of $\tau_*$ is the same as the datum of a natural
transformation $\tau^*:g^*\Rightarrow f^*$.

(iv)\ \
Let $\sU'$ be a universe $\sU\in\sU'$; the $(\sU,\sU')$-topoi
are the objects of a $2$-category, denoted
$$
(\sU,\sU')\tdu\Topos
$$
whose $1$-cells are the morphisms of topoi, and whose
$2$-cells are the natural transformations between such
morphisms. Obviously this category depends on the choice
of $\sU'$, but two different choices are related by a
$2$-equivalence of $2$-categories, so we shall usually
omit mentioning an explicit choice of $\sU'$, and write
$\sU\tdu\Topos$ for this $2$-category. Moreover, when
the choice of $\sU$ is clear from the context, we shall
likewise drop the mention of $\sU$ and write just $\Topos$.
\end{definition}

\begin{example}\label{ex_coproduct-of-topoi}
Let $T_\bullet:=(T_i~|~i\in I)$ be any $\sU$-small family of
$(\sU,\sU')$-topoi.

(i)\ \
The product category $T:=\prod_{i\in I}T_i$ is a topos
(see example \ref{ex_cat-cats}(i)). Indeed, pick for
each $i\in I$ a small site $C_i:=(\cC_i,J_i)$ with an
equivalence $C_i^\sim\isom T_i$; according to example
\ref{ex_cat-cats}(i), the coproduct of the family of
categories $(\cC_i~|~i\in I)$, is represented by a
category $\cC$ whose set of objects is the disjoint
union $\amalg_{i\in I}\Ob(\cC_i)$. We endow $\cC$ with
the topology $J$ such that $J(i,X):=J_i(X)$ for every
$(i,X)\in\Ob(\cC)$. Then it is easily seen that there
is a natural isomorphism of categories
$$
(\cC,J)^\sim\isom\prod_{i\in I}C_i^\sim.
$$
However, the product of the categories $C_i^\sim$ also
represents the $2$-product, in the $2$-category $\sU'\tdu\bCat$,
of the same family of categories, and the standard universal
cone is also a universal pseudo-cone (see example
\ref{ex_2-products}); there follows an equivalence
of categories
$$
(\cC,J)^\sim\isom T
$$
whence the contention.

(ii)\ \
Moreover, $T$ represents the $2$-coproduct of the family
$T_\bullet$ in the $2$-category $\Topos$. Indeed,
for every $i\in I$ let $\pi_i:T\to T_i$ be the projection
functor. We define a right adjoint $e_i:T_i\to T$ for
$\pi_i$ by the rule : $X\mapsto(X_j~|~j\in I)$ for every
$X\in\Ob(T_i)$, where $X_i:=X$ and $X_j:=1_{T_j}$ is a
fixed choice of a final object in $T_j$, for every
$j\in I\setminus\{i\}$. It is easily seen that $\pi_i$
is left exact, and there is a natural choice of adjunction
for the pair $(\pi_i,e_i)$, whence a well defined morphism
of topoi
$$
\phi_i:T_i\to T
\qquad
\text{such that $\phi_i^*=\pi_i$ and $\phi_{i*}:=e_i$}.
$$
Next, let $S$ be another topos, and $(\psi_i:T_i\to S~|~i\in I)$
a given family of morphisms of topoi. By the universal property
of the product, there exists a unique functor $\lambda^*:S\to T$
such that
\set\begin{equation}\label{eq_topoi-dolorosi}
\phi^*_i\circ\lambda^*=\psi^*_i
\qquad
\text{for every $i\in I$}
\end{equation}
and it is easily seen that $\lambda^*$ is left exact. We
construct a right adjoint $\lambda_*:T\to S$ for $\lambda^*$
by the following rule. Given any object
$X_\bullet:=(X_i~|~i\in I)$ of $T$, let $\lambda_*(X_\bullet)$
be any object of $S$ representing the product
$\prod_{i\in I}\psi_{i*}X_i$, and fix as well a universal cone
$(p^X_i:\lambda_*(X_\bullet)\to\psi_{i*}X_i~|~i\in I)$.
Then, for every morphism $f_\bullet:X_\bullet\to Y_\bullet$ in
$T$ there is a unique morphism
$\lambda_*(f_\bullet):\lambda_*(X_\bullet)\to\lambda_*(Y_\bullet)$
that makes commute the diagram
$$
{\diagram
\lambda_*(X_\bullet) \ar[rr]^-{\lambda_*(f_\bullet)} \ar[d]_{p^X_i}
& & \lambda_*(Y_\bullet) \ar[d]^{p^Y_i} \\
\psi_{i*}X_i \ar[rr]^-{\psi_{i*}f_i} & & \psi_{i*}Y_i
\enddiagram}
\qquad
\text{for every $i\in I$}.
$$
Thus, the rules $X_\bullet\mapsto\lambda_*(X_\bullet)$ and
$f_\bullet\mapsto\lambda_*(f_\bullet)$ yield a well defined
functor as sought, unique up to unique isomorphism.
Then it is easily seen that $\lambda^*$ is left adjoint
to $\lambda_*$, and there is even a natural choice of
adjunction for the pair $(\lambda^*,\lambda_*)$ (details
left to the reader), so the latter yields a well defined
morphism of topoi $\lambda:T\to S$. Lastly,
\eqref{eq_topoi-dolorosi} implies that $\lambda_*\circ\phi_*$
is right adjoint to $\psi^*_i$, whence an isomorphism of
functors
$$
\lambda_*\circ\phi_*\isom\psi_{i*}
\qquad
\text{for every $i\in I$}
$$
which shows that the system $(\phi_i~|~i\in I)$ is a
universal pseudo-cocone, whence the contention.
\end{example}

\begin{proposition}\label{prop_half-dominates}
Let $T$, $T'$ be two topoi, $f:T\to T'$ a left exact functor,
$\phi$ a morphism of\/ $T$.
\begin{enumerate}
\item
If $f$ commutes with all small colimits in $T$, the following
holds :
\begin{enumerate}
\item
$f$ is continuous for the canonical topologies on $T$ and $T'$.
\item
There exists a morphism of topoi $F:T'\to T$, unique up to
unique isomorphism, such that $F^*=f$.
\end{enumerate}
\item
$f$ is conservative if and only if it reflects epimorphisms.
\item
If $f$ is exact, the natural morphism
$f(\Img\,\phi)\to\Img(f\phi)$ is an isomorphism.
\end{enumerate}
\end{proposition}
\begin{proof}(i.a): Let $S:=\{g_i:X_i\to X~|~i\in I\}$ be a small
covering family of morphisms in $\Can_T$; by lemma
\ref{lem_crit-continuity}, it suffices to show that
$fS:=\{f(g_i)~|~i\in I\}$ is a covering family in $\Can_{T'}$.
However, our assumption implies that
the induced morphism $\amalg_{i\in I}X_i\to X$ is an
epimorphism, hence the same holds for the induced morphism
$\amalg_{i\in I}fX_i\to fX$ in $T'$ (by proposition
\ref{prop_was-also-cofinal}(ii)). But all epimorphisms are
universal effective in $T'$ (remark \ref{rem_summarized}(i)),
hence $(f(g_i)~|~i\in I)$ is a universal effective family, as
required.

(i.b): By corollary \ref{cor_two-U-sites}(i.a,i.c), the functor
$f$ gives rise to a morphism of topoi
$$
(\tilde f{}^*,\tilde f_*,\eta):\Can(T')^\sim\to\Can(T)^\sim
$$
(where $\eta$ is any choice of unit of adjunction), and it
suffices to show that $\tilde f{}^*$ is isomorphic to
$f$, under the natural identifications $T\isom\Can(T)^\sim$,
$T'\isom\Can(T')^\sim$ given by the Yoneda embeddings.
The latter follows from corollary \ref{cor_two-U-sites}(i.b).

(iii) follows easily from remark \ref{rem_mono-epi-fact}.

(ii) follows from propositions \ref{prop_univ-effective}(i)
and \ref{prop_conser-nd-left-exact}.
\end{proof}

\begin{remark}\label{rem_Topos-*}
(i)\ \
Let us consider the $2$-category $(\sU,\sU')\tdu\Topos^*$ whose
objects are the same as those of the $2$-category
$(\sU,\sU')\tdu\Topos$, and such that for every two $\sU$-topoi
$T$ and $S$, the $1$-cells $f:T\to S$ are the left exact functors
$f:S\to T$ that commute with all small colimits. For any pair of
such functors $f,f':S\to T$, the $2$-cells $\beta:f\Rightarrow f'$
of $(\sU,\sU')\tdu\Topos^*$ are the natural transformations
$\beta:f'\Rightarrow f$. As usual, we shall often write
$\sU\tdu\Topos^*$ or just $\Topos^*$ for this $2$-category.
We have a natural strict pseudo-functor
$$
(-)^*:\Topos\to\Topos^*
$$
that is the identity on objects, and assigns to every
morphism $f:=(f^*,f_*,\eta_f):T\to S$ of topoi the functor
$f^*:S\to T$, and to every natural transformation
$\beta:(f^*,f_*,\eta_f)\Rightarrow(g^*,g_*,\eta_g)$ the adjoint
transformation $(\beta,\eta_f,\eta_g)^\dagger:g^*\Rightarrow f^*$.
Then proposition \ref{prop_half-dominates}(i.b) can be restated
by saying that $(-)^*$ is a strict $2$-equivalence of
$2$-categories.

(ii)\ \
Notice also the strict pseudo-functors
$$
\sU'\tdu\bCat^o\xleftarrow{(-)^*}\Topos
\xrightarrow{(-)_*}\sU'\tdu\bCat
$$
defined as follows. The pseudo-functor $(-)_*$ is the identity
map on objects : $T\mapsto T$, and is given on morphisms by
the rule : $(F_*,F^*,\eta_F)\mapsto F_*$ (notation of definition
\ref{def_morph-of-topoi}(i)). On natural transformations it
is also the identity : $(\tau:(F^*,F_*,\eta_F)\Rightarrow
(G^*,G_*,\eta_G))\mapsto(\tau:F_*\Rightarrow G_*)$.

The pseudo-functor $(-)^*$ is given on objects by the rule :
$T\mapsto T^o$, and on morphisms by the rule :
$(F^*,F_*,\eta_F)\mapsto(F^*)^o$. On natural transformations
it is defined by the rule :
$$
(\tau:(F^*,F_*,\eta_F)\Rightarrow(G^*,G_*,\eta_G))\mapsto
((\tau,\theta_F,\theta_G)^{\dagger o}:F^{*o}\Rightarrow G^{*o})
$$
where $\theta^F$ and $\theta^G$ are the adjunctions for the
pairs $(F^*,F_*)$ and respectively $(G^*,G_*)$ whose units
are $\eta_F$ and respectively $\eta_G$, and
$(\tau,\theta_F,\theta_G)^\dagger$ is the adjoint to $\tau$,
as in remark \ref{rem_adjoint-transf}(ii). The pseudo-functoriality
of $(-)^*$ follows straightforwardly from remark
\ref{rem_adjoint-transf}(iv,v) : the details shall be left
to the reader.
\end{remark}

\begin{example}\label{ex_another-basic}
Let $T$ be a topos.

(i)\ \
Say that $T=(\cC,J)^\sim$ for a small site $(\cC,J)$.
According to remark \ref{rem_sheaves}(vi), the category
$\cC^\wedge$ is also a topos, and since the functor
$F\mapsto F^a$ of theorem \ref{th_ass-sheaf} is left exact,
it determines a morphism of topoi $T\to\cC^\wedge$. (On the
other hand, the category $T^\wedge$ is too large to be a
$\sU$-topos.)

(ii)\ \
Let $f:=(f^*,f_*,\eta):T\to S$ be any morphism of topoi;
we have an essentially commutative diagram whose horizontal
arrows are the Yoneda embeddings :
$$
\xymatrix@C+20pt{
T \ar[r]^-{h_T} \ar[d]_{f_*} &
\Can(T)^\sim \ar[d]^{(f^*)^\sim_*} \\
S \ar[r]^-{h_S} & \Can(S)^\sim.
}$$
Indeed, an isomorphism $(f^*)^\sim_*\circ h_T\isom h_S\circ f_*$
is given explicitly as follows. Let $\theta$ be the adjunction
for the pair $(f^*,f_*)$ resulting from the unit $\eta$. So,
$\theta$ consists of a system of bijections:
$$
h_X\circ f^*(Y)=\Hom_T(f^*Y,X)\xrightarrow{\theta_{Y,X}}
\Hom_S(Y,f_*X)=h_{f_*X}(Y)
$$
natural in $X\in\Ob(T)$ and $Y\in\Ob(S)$. The naturality
with respect to morphisms $Y'\to Y$ in $S$ then implies
that the rule $Y\mapsto\theta_{Y,X}$ is an isomorphism of
sheaves $\theta_{\bullet,X}:h_X\circ f^*\isom h_{f_*X}$ for
every $X\in\Ob(T)$. Lastly, the naturality with respect
to morphisms $X\to X'$ in $T$ says that the rule
$X\mapsto\theta_{\bullet,X}$ yields an isomorphism of functors
as sought : the details are left to the reader.
\end{example}

\sset\subsubsection{Global section functors}
\label{subsec_glob-sections}
In view of theorem \ref{th_canon-topos}(iv), the objects of
any topos $T$ may be thought of as sheaves on $\Can(T)$, and
one uses often the suggestive notation
$$
X(S):=\Hom_T(S,X) \qquad \text{for every $X,S\in\Ob(T)$.}
$$
The elements of $X(S)$ are also called the {\em $S$-sections\/} of
$X$. If $1_T$ is any final object of $T$, one defines the
{\em global sections functor\/} $\Gamma:T\to\Set$ by the rule :
$$
U\mapsto\Gamma(T,U):=U(1_T).
$$
Moreover $\Gamma$ admits a left adjoint :
$$
(-)_T:\Set\to T\ :\ S\mapsto S_T:=S\times 1_T
$$
(the coproduct of $S$ copies of $1_T$) and one calls $S_T$ the
{\em constant sheaf with value $S$}.

\begin{proposition}\label{prop_glob-sections}
With the notation of \eqref{subsec_glob-sections}, we have :

{\em(i)}\ \
The functor $(-)_T$ is exact, hence the pair $((-)_T,\Gamma)$
yields a morphism of topoi
$$
\Gamma_T:T\to\Set.
$$

{\em(ii)}\ \
For every pair of morphisms of topoi $u,v:T\to\Set$ there
exists a unique natural transformation $\beta:u_*\Rightarrow v_*$,
and $\beta$ is an isomorphism of functors.

{\em(iii)}\ \
For every morphism $f:T\to S$ of topoi, there exists a unique
isomorphism :
$$
\Gamma_S\circ f\isom\Gamma_T.
$$
\end{proposition}
\begin{proof} (See \cite[Exp.IV, \S4.3]{SGA4-1}). For the proof
of (i), we may assume that $T=C^\sim$ for a small site $C:=(\cC,J)$.
Then, for every (small) set $S$, let $c_S:\cC^o\to\Set$ be the
constant presheaf on $\cC$ with value $S$; since the colimits in
$\cC^o$ are computed argumentwise (corollary \ref{cor_pre-misc}(ii)),
we have $c_S=\coprod_{s\in S}c_{\{s\}}$, whence $S_T=c_S^a$, since the
functor $(-)^a$ commutes with small colimits. Since $(-)^a$
is an exact functor (theorem \ref{th_ass-sheaf}(ii)), it suffices
to check that the same holds for the functor
$\Set\to\cC^o$ : $S\mapsto c_S$. But more precisely, the latter
commutes with all small limits and all small colimits, again
because the latter are computed argumentwise.

(ii): Since $u^*$ is exact, $u^*(\{\emptyset\})$ is a final
object of $T$, and we get natural bijections :
$$
u_*X\isom\Hom_\Set(\{\emptyset\},u_*X)\isom
\Hom_T(u^*(\{\emptyset\}),X)
\qquad
\text{for every $X\in\Ob(T)$}
$$
{\em i.e.} the presheaf $u_*:T\to\Set$ on the category $T^o$
is represented by the initial object $u^*(\{\emptyset\})^o$.
The same applies to $v_*$, and then Yoneda's lemma yields
natural bijections
$$
\Hom_{(T^o)^\wedge}(u_*,v_*)\isom
\Hom_T(v^*(\{\emptyset\}),u^*(\{\emptyset\}))\isom\{\emptyset\}
$$
whence the contention.

(iii) is an immediate consequence of (ii).
\end{proof}

\begin{remark}\label{rem_choose-two-univs}
(i)\ \
If $T$ is a $\sU$-topos, $C:=(\cC,J)$ is a $\sU$-lex-site,
and $g:\Can(T)\to C$ is a morphism of sites, then the
underlying functor $g:\cC\to T$ is left exact. Indeed,
by corollary \ref{cor_two-U-sites}(i.b) we have an
essentially commutative diagram of functors whose vertical
arrows are the Yoneda embeddings :
\set\begin{equation}\label{eq_sites-to-topoi}
{\diagram \cC \ar[r]^-g \ar[d]_{h_\cC^a} & T \ar[d]^{h_T} \\
C^\sim \ar[r]^-{\tilde g{}^*} & \Can(T)^\sim.
\enddiagram}
\end{equation}
Then, $h_T$ is an equivalence (theorem \ref{th_canon-topos}(iv))
and $\tilde g{}^*$ is left exact, since $g$ is a morphism
of sites; it suffices therefore to show that $h_\cC^a$
is left exact. But the functor $(-)^a:\cC^\wedge\to C^\sim$
is left exact by theorem \ref{th_ass-sheaf}(ii), so the
assertion follows from corollary \ref{cor_pre-misc}(vi).

(ii)\ \
Every morphism of $\sU$-sites $u:C'\to C$ induces a morphism
of topoi
$$
\tilde u:=(\tilde u{}^*,\tilde u_*,\tilde\eta{}^u):
C^{\prime\sim}\to C^\sim
$$
after choosing a unit $\tilde\eta{}^u$ for the adjoint pair
$(\tilde u{}^*,\tilde u_*)$ (corollary \ref{cor_two-U-sites}(i)).

(iii)\ \
Conversely, a morphism $(f^*,f_*,\eta):T\to S$ of topoi
determines a morphism of sites
$$
f^*:\Can(T)\to\Can(S).
$$
Indeed, $g:=f^*:S\to T$ is continuous for the canonical
topologies, by virtue of proposition \ref{prop_half-dominates}(i),
and by corollary \ref{cor_two-U-sites}(i.b) we have an
essentially commutative diagram of functors as in
\eqref{eq_sites-to-topoi}, with $\cC:=S$ and $C=\Can(S)$.
In this case, both vertical arrows of \eqref{eq_sites-to-topoi}
are equivalences (theorem \ref{th_canon-topos}(iv)). Since
$g$ is left exact, it then follows that the same holds
for $\tilde g{}^*$, whence the assertion.

(iv)\ \
Likewise, for every pair of $\sU$-sites $C:=(\cC,J),C':=(\cC',J')$,
for every cocontinuous functor $v:\cC\to\cC'$, it is clear that
$\breve v^*$ is an exact functor, so $v$ induces a morphism of topoi
$$
\breve v:=(\breve v{}^*,\breve v_*,\breve\eta{}^v):C^\sim\to C'^\sim
$$
(corollary \ref{cor_two-U-sites}(ii)), after choosing a unit
$\breve\eta{}^v$ for the adjoint pair $(\breve v^*,\breve v_*)$.

(v)\ \
In the situation of (iv), suppose that $v$ admits a right adjoint
$u$. Then $u$ is continuous for the topologies of the sites $C$
and $C'$, and $\tilde u{}^*$ is isomorphic to $\breve v{}^*$
(lemma \ref{lem_needed}). Especially, $\tilde u{}^*$ is exact,
so $u$ is a morphism of sites, and we have an isomorphism of
morphisms of topoi :
$$
\tilde u\isom\breve v.
$$
\end{remark}

\sset\subsubsection{}\label{subsec_T-and-Can}
We define as follows a pair of pseudo-functors
$$
\xymatrix{ \Site \ar@<.5ex>[rr]^-\sT & &
\Topos. \ar@<.5ex>[ll]^-{\Can}
}$$
Let $C$ be a $\sU$-site; by remark \ref{rem_summarized}(iv)
we may find an isomorphism of topoi
$$
\omega_C:C^\sim\isom\sT(C)
$$
with a $\sU$-topos $\sT(C)$, and it is easily seen that
we may take this topos to be also $\sU'$-small. Then the rule :
$C\mapsto\sT(C)$ defines the pseudo-functor $\sT$ on objects.
Next, let $g:C\to C'$ be any morphism of $\sU$-sites; by remark
\ref{rem_choose-two-univs}(ii), any choice of a unit $\eta^g$
for the adjoint pair $(\tilde g{}^*,\tilde g_*)$ yields a
morphism of topoi $\tilde g:C^\sim\to C'^\sim$.
There follows easily a unique morphism of $\sU$-topoi
$$
\sT(g):=(\sT(g)^*,\sT(g)_*,\eta^{\sT(g)}):\sT(C)\to\sT(C')
$$
that is identified with $(\tilde g{}^*,\tilde g_*,\eta^g)$
via the isomorphisms $\omega_C$ and $\omega_{C'}$
(details left to the reader). Lastly, let
$\beta:g\Rightarrow g'$ be a natural transformation between
morphisms of sites $g,g':C\to C'$; there follows
a natural transformation $\beta^\wedge:g'^\wedge\Rightarrow g^\wedge$
(notation of \eqref{eq_pullback-presheaves}) that induces by
restriction a natural transformation
$$
\beta^\sim_*:\tilde g'_*\Rightarrow\tilde g_*
$$
which is in turn identified via $\omega_C$ and
$\omega_{C'}$ with a unique natural transformation
$$
\sT(\beta):\sT(g')_*\Rightarrow\sT(g)_*.
$$
If $g'':C\to C'$ is a third morphism of sites and
$\beta':g'\Rightarrow g''$ is another natural transformation,
it is clear that
$$
\sT(\beta'\odot\beta)=\sT(\beta)\odot\sT(\beta').
$$
Likewise, if $h,h':C'\to C''$ is another pair of morphisms
of sites and $\alpha:h\Rightarrow h'$ is another natural
transformation, we deduce easily from
\eqref{eq_pullback-presheaves} the identity :
$$
\sT(\alpha*\beta)=\sT(\beta)*\sT(\alpha).
$$
For a composable pair
$C\xrightarrow{g}C'\xrightarrow{g'}C''$ of morphisms of
$\sU$-sites, it is easily seen that
$\sT(g')_*\circ\sT(g)_*=\sT(g'\circ g)_*$, and the coherence constraint
$\gamma^\sT_{g,g'}$ of $\sT$ is given by $\one_{\sT(g'\circ g)_*}$. Likewise,
$\sT(\one_C)_*=\one_{\sT(C)}$ for every $\sU$-site $C$,
and the coherence constraint $\delta^\sT_C$ is given by the
identity automorphism of $\one_{\sT(C)}$. Then the required
coherence axioms are trivially satified; however, notice that
$\sT$ is {\em not} a strict pseudo-functor, since we do not have
necessarily $\sT(g)^*\circ\sT(g')^*=\sT(g'\circ g)^*$. This
completes the construction of $\sT$.

The pseudo-functor $\Can$ assigns to every $\sU$-topos $T$
the $\sU$-site $\Can(T)$ and to every morphism of $\sU$-topoi
$f:T\to S$ the morphism of sites $f^*:\Can(T)\to\Can(S)$
(see remark \ref{rem_choose-two-univs}(iii)). Lastly, let
$(f^*,f_*,\eta^f),(g^*,g_*,\eta^g):T\to S$ be two morphisms of
topoi and $\beta:f_*\Rightarrow g_*$ a natural transformation;
then we let
$\Can(\beta):=(\beta,\eta^f,\eta^g)^\dagger:g^*\Rightarrow f^*$,
the adjoint transformation of $\beta$, as defined in remark
\ref{rem_adjoint-transf}(ii). Taking into account remark
\ref{rem_adjoint-transf}(iv,v,vi), we easily see that these
rules define a strict pseudo-functor $\Can$ as sought.

\begin{remark}\label{rem_lex-variants}
Notice that the pseudo-functor $\Can$ factors through the
inclusion strict pseudo-functor $\lex.\Site\to\Site$. We
then get as well a pair of pseudo-functors :
$$
\xymatrix{ \lex.\Site \ar@<.5ex>[rr]^-{\lex.\sT} & &
\Topos. \ar@<.5ex>[ll]^-{\lex.\Can}
}$$
where $\lex.\sT$ is the restriction of $\sT$.
\end{remark}

\begin{theorem}\label{th_adj-topos-site}
{\em(i)}\ \
The pseudo-functor $\sT$ is right $2$-adjoint to the
pseudo-functor $\Can$.

{\em(ii)}\ \
The pseudo-functor $\lex.\sT$ is right $2$-adjoint to the
pseudo-functor $\lex.\Can$.

{\em(iii)}\ \
The pseudo-functors $\Can$ and $\lex.\Can$ are fully faithful.
\end{theorem}
\begin{proof}(i): We construct as follows a pseudo-natural
equivalence $\shh_\bullet:\one_\Topos\Rightarrow\sT\circ\Can$.
For every $\sU$-topos $T$, we choose a quasi-inverse
$\shh^*_T:\sT\circ\Can(T)\to T$ for the composition
$\shh_{T*}:T\to\sT\circ\Can(T)$ of the Yoneda embedding
$h_T:T\to\Can(T)^\sim$ with the isomorphism $\omega_{\Can(T)}$
of \eqref{subsec_T-and-Can}. We also fix an isomorphism of
functors
$\eta^{\shh_T}:\one_{\sT\circ\Can(T)}\isom\shh_{T*}\circ\shh^*_T$,
and we consider the morphism of $\sU$-topoi :
$$
\shh_T:=(\shh_T^*,\shh_{T*},\eta^{\shh_T}):T\to\sT\circ\Can(T).
$$
In order to show that the rule $T\mapsto\shh_T$ defines the
sought pseudo-natural equivalence, we need to explicit its
coherence constraint; the latter amounts to an isomorphism
of functors
$$
\tau^\shh_f:\sT(f^*)_*\circ\shh_{T*}\isom\shh_{S*}\circ f_*
$$
for every morphism of $\sU$-topoi $(f^*,f_*,\eta^f):T\to S$.
However, example \ref{ex_another-basic}(ii) yields an isomorphism
of functors $\tau^h_f:(f^*)^\sim_*\circ h_T\isom h_S\circ f_*$,
so we can take :
$$
\tau^\shh_f:=\omega_{(S,\Can_S)}*\tau^h_f.
$$
To check the coherence axioms for $\tau^\shh$, notice first
that $\sT$ and $\Can$ are unital pseudo-functors, so that
the first coherence axiom amounts to the identity :
$\tau^\shh_{\one_T}=\one_{\shh_{T*}}$ for every $\sU$-topos $T$.
The latter follows from the identity : $\tau^h_{\one_T}=\one_{h_T}$,
which is clear from the construction of $\tau^h$.

Next, for two morphisms of $\sU$-topoi $f:=(f^*,f_*,\eta^f):T\to T'$,
$g:=(g^*,g_*,\eta^g):T'\to T''$, we come down to checking that
\set\begin{equation}\label{eq_combine-adjunctions}
(\eta^h_g*f_*)\odot((g^*)^\sim_**\eta^h_f)=\eta^h_{g\circ f}.
\end{equation}
Now, recall that $\eta^h_f$ is induced by the unique adjunction
$\theta^f$ for the pair $(f^*,f_*)$ whose unit is $\eta^f$, and
likewise for $\eta^h_g$ and $\eta^h_{g\circ f}$; on the other hand,
$\eta^h_{g\circ f}$ is precisely the unit of the adjunction
$\theta^g\circ\theta^f$ (notation of remark
\ref{rem_adjoint-transf}(i)). From this, the identity
\eqref{eq_combine-adjunctions} follows straightforwardly.

Likewise, let $f,f':T\to S$ be two morphisms of topoi, and
$\beta:f\Rightarrow f'$ a natural transformation; in order
to check the naturality of $\tau^\shh$ we come down to showing
that :
$$
\tau^h_{f'}\odot(\beta^\dagger)^\sim_*=(h_S*\beta)\odot\tau^h_f.
$$
The latter follows by inspecting the characterization of
$\beta^\dagger$ furnished by remark \ref{rem_adjoint-transf}(iii),
and this completes the construction of $\shh_\bullet$. From the
pseudo-natural equivalence $\shh_\bullet$ we construct the sought
$2$-adjunction as the composition of the pseudo-natural
transformation
$$
H_\sT*(\Can^o\times\one_\Site):
\Hom_\Site(\Can,-)\Rightarrow\Hom_\Topos(\sT\circ\Can,\sT)
$$
with the pseudo-natural transformation
$$
H_\Topos*(\shh^o_\bullet\times\one_\Topos):
\Hom_\Topos(\sT\circ\Can,\sT)\Rightarrow\Hom_\Topos(\one_\Topos,\sT).
$$
This composition assigns to every $\sU$-topos $T$ and every
$\sU$-site $C$ a functor
$$
\Phi_{T,C}:\Hom_\Site(\Can(T),C)\to\Hom_\Topos(T,\sT(C))
$$
as required, and it remains to check that this functor is
an equivalence for every such $T$ and $C$.

Explicitly, $\Phi_{T,C}$ assigns to every morphism of sites
$g:\Can(T)\to C$ the morphism of topoi
$\sT(g)\circ\shh_T:T\to\sT(C)$, and to every $2$-cell
$\beta:g\Rightarrow g'$ of $\Site$ the natural transformation
$\sT(\beta)*\shh_T:
\sT(g)\circ\shh_T\Rightarrow\sT(g')\circ\shh_T$.
Let us show the essential surjectivity : thus, say that
$C=(\cC,J)$, and let $f:=(f^*,f_*,\eta^f):T\to\sT(C)$ be
a morphism of topoi; from theorem \ref{th_canon-topos}(ii)
we see that $g:=f^*\circ\omega_C\circ h^a_\cC:\cC\to T$ is a
morphism of sites $g:\Can(T)\to C$, and we have an isomorphism
of functors $h_T\circ g\isom\tilde g{}^*\circ h^a_\cC$
(corollary \ref{cor_two-U-sites}(i.b)). We deduce an
isomorphism :
\set\begin{equation}\label{eq_iso-with-h^a}
\shh_{T*}\circ f^*\circ\omega_C\circ h^a_\cC\isom
\omega_{\Can(T)}\circ\tilde g{}^*\circ h^a_\cC.
\end{equation}

\begin{claim}\label{cl_remove-h^a}
Let $D$ be any category, $k,k':C^\sim\to D$
two functors, $\lambda:k\circ h^a_\cC\isom k'\circ h^a_\cC$
an isomorphism of functors, and suppose that $k$ commutes
with all small colimits. Then there exists an isomorphism
of functors $\lambda':k\isom k'$.
\end{claim}
\begin{pfclaim} Let $G\subset\Ob(\cC)$ be a small topologically
generating subset, $\cG\subset\cC$ the full subcategory with
$\Ob(\cG)=G$; endow $\cG$ with the topology $J_\cG$ induced by
$J$ via the inclusion functor $u:\cG\to\cC$, and set
$G:=(\cG,J_\cG)$. Then $\tilde u_*:C^\sim\to G^\sim$ is an
equivalence (proposition \ref{prop_comparis-lemma}), hence the
same holds for its left adjoint $\tilde u{}^*:G^\sim\to C^\sim$.
According to corollary \ref{cor_two-U-sites}(i.b), we have an
isomorphism of functors
$\mu:\tilde u{}^*\circ h^a_\cG\isom h^a_\cC\circ i$. There follows
a unique isomorphism $\lambda':k\circ\tilde i{}^*\circ h^a_\cG
\isom k'\circ\tilde i{}^*\circ h^a_\cG$ such that
$$
(k'*\mu)\odot\lambda'=(\lambda*i)\odot(k*\omega).
$$
Now, if $k\circ\tilde i{}^*$ and $k'\circ\tilde i{}^*$ are
isomorphic, the same follows for $k$ and $k'$, since
$\tilde i{}^*$ is an equivalence (details left to the reader).
Thus, we may replace $C$ by $G$ and $\lambda$ by $\lambda'$,
and assume from start that $C$ is a small site.
Next, let $F\in\Ob(C^\sim)$; we have a universal cocone
$\tau^F$ given by the rule :
$$
(X,s)\mapsto(\tau^F_{(X,s)}:h^a_X\to F)
\qquad
\text{for every $(X,s)\in\Ob(\cFib(F))$}
$$
where $\tau^F_{(X,s)}$ is characterized as the morphism of
sheaves whose composition with the natural morphism of
presheaves $h_X\to h^a_X$ is the unique morphism of presheaves
$t^F_{(X,s)}:h_X\to F$ such that $(t^F_{(X,s)})_X(\one_X)=s$
(remark \ref{rem_rep-and-sheafify}(iii)). By assumption,
the cocone $k*\tau^F$ is still universal, hence there exists
a unique isomorphism $\lambda'_F:kF\isom k'F$ that makes
commute the diagram :
$$
{\spreaddiagramcolumns{+20pt}\diagram
k(h^a_X) \ar[r]^-{k(\tau^F_{(X,s)})} \ar[d]_{\lambda_X} &
kF \ar[d]^{\lambda'_F} \\
k'(h^a_X) \ar[r]^-{k'(\tau^F_{(X,s)})} & k'F
\enddiagram}
\qquad
\text{for every $(X,s)\in\Ob(\cFib(F))$}
$$
and we are reduced to checking that the rule : $F\mapsto\lambda'_F$
is a natural transformation $k\Rightarrow k'$. Thus, let
$\phi:F\to G$ be any morphism of sheaves on $C$; recall that
$\phi$ induces a functor
$$
\cFib(\phi):\cFib(F)\to\cFib(G)
\qquad
(X,s)\mapsto(X,\phi_X(s))
$$
and the foregoing characterizations of $\tau^F$ and $\tau^G$
easily imply that :
\set\begin{equation}\label{eq_identity}
\phi\circ\tau^F_{(X,s)}=\tau^G_{(X,\phi_X(s))}
\qquad
\text{for every $(X,s)\in\Ob(\cFib(F))$}.
\end{equation}
Consider now for every $(X,s)\in\Ob(\cFib(F))$ the diagram :
$$
\xymatrix@C+20pt{
k'(h^a_X) \ar[rrr]^-{k'(\tau^F_{(X,s)})} \ar@{=}[ddd] & & &
k'F \ar[ddd]^{k'(\phi)} \\
& k(h^a_X) \ar[r]^-{k(\tau^F_{(X,s)})} \ddouble \ar[lu]_{\lambda_X} &
kF \ar[d]^{k(\phi)} \ar[ru]^{\lambda'_F} \\
& k(h^a_X) \ar[r]^-{k(\tau^G_{(X,\phi(s))})} \ar[ld]_{\lambda_X} &
kG \ar[rd]^{\lambda'_G} \\
k'(h^a_X) \ar[rrr]^-{k'(\tau^G_{(X,\phi(s))})} & & & k'G
}$$
We have to prove the commutativity of the right trapezoidal
subdiagram, and by \eqref{eq_identity}, the inner and outer
square subdiagrams both commute; moreover, also the upper
and lower trapezoidal subdiagrams commute, by construction.
Now, since by assumption the cocone $k*\tau^F$ is still
universal, it suffices to check that
$k'(\phi)\circ\lambda'_F\circ k(\tau^F_{(X,s)})=
\lambda'_G\circ k(\phi)\circ k(\tau^F_{(X,s)})$. The latter
follows from a straightforward diagram chase, left to the
reader.
\end{pfclaim}

Since $\tilde g{}^*$ commutes with small colimits, claim
\ref{cl_remove-h^a} and \eqref{eq_iso-with-h^a} yield an
isomorphism :
$$
\shh_{T*}\circ f^*\circ\omega_C\isom\omega_{\Can(T)}\circ\tilde g{}^*.
$$
From this, we further deduce an isomorphism of functors :
$$
f^*\isom\shh^*_T\circ\sT(g)^*
$$
whence, finally, an isomorphism of their respective right
adjoint functors : $\sT(g)_*\circ\shh_{T*}\isom f_*$.

Next, we check that $\Phi_{T,C}$ is faithful : let
$g,g':\Can(T)\to C$ be two morphisms of sites, and
$\beta,\beta':g\Rightarrow g'$ two $2$-cells of $\Site$
such that $\sT(\beta)_**\shh_T=\sT(\beta')_**\shh_T$; we
need to show that $\beta=\beta'$. However, the condition
means that $\tilde\beta_*\circ h_T=\tilde\beta{}'_*\circ h_T$;
{\em i.e.} for every $X\in\Ob(T)$, every $A\in\Ob(\cC)$
and every morphism $\phi:gA\to X$ in $T$, we have
$\phi\circ\beta_A=\phi\circ\beta'_A$. Letting $X:=gA$
and $\phi:=\one_{gA}$, we get $\beta_A=\beta'_A$, for
every $A\in\Ob(\cC)$, as required.

Lastly, in order to check that $\Phi_{T,C}$ is full, let
$\alpha:\sT(g)\circ\shh_T\Rightarrow\sT(g')\circ\shh_T$ be
any natural transformation, with $g$ and $g'$ as in the
foregoing. This is the same as a natural transformation
$\alpha':\tilde g_*\circ h_T\Rightarrow\tilde g{}'_*\circ h_T$.
The latter assigns to every $X\in\Ob(T)$ a morphism of
sheaves $\alpha'_X:h_X\circ g^o\to h_X\circ g'^o$; in
particular, for every $A\in\Ob(\cC)$, we deduce a map
$$
(\alpha'_{gA})_A:\Hom_T(gA,gA)\to\Hom_T(g'A,gA)
\qquad\text{and we set :}\qquad
\beta_A:=(\alpha'_{gA})_A(\one_{gA}).
$$
Let us check that the rule : $A\mapsto\beta_A$ yields a
natural transformation $\beta:g'\to g$. Indeed, let
$\phi:A\to B$ be any morphism of $\cC$; we need to show
that $g(\phi)\circ\beta_A=\beta_B\circ g'(\phi)$. However,
the naturality of the rule $X\mapsto\alpha'_X$ for every
$X\in\Ob(T)$ implies the identity :
$$
g(\phi)\circ(\alpha'_{gA})_A(\one_{gA})=
(\alpha'_{gB})_A(g(\phi))
$$
and the naturality of the rule : $Y\mapsto(\alpha'_{gB})_Y$
for every $Y\in\Ob(\cC)$ yields the identity :
$$
(\alpha'_{gB})_A(g(\phi))=(\alpha'_{gB})_B(\one_{gB})\circ g'(\phi)
$$
whence the contention. It remains to check that
$\sT(\beta)_**\shh_T=\alpha$, or equivalently, that
$\tilde\beta_**h_T=\alpha'$. The latter comes down to the
identity : $(\alpha'_X)_A(\phi)=\phi\circ\beta_A$ for every
$X\in\Ob(T)$, every $A\in\Ob(\cC)$ and every morphism
$\phi:gA\to X$ in $T$. This identity in turn follows easily
from the naturality of the rule : $X\mapsto\alpha'_X$ (details
left to the reader).

(ii): Notice that -- by virtue of example
\ref{ex_simple-case} -- the same rules defining
$\shh_\bullet$ also yield a pseudo-natural transformation
$\one_\Topos\Rightarrow\lex.\sT\circ\lex.\Can$. We use
this pseudo-natural transformation as in the foregoing to
construct the sought $2$-adjunction for the pair
$(\lex.\Can,\lex.\sT)$; the details shall be left to the
reader.

(iii): Since $\shh_\bullet$ is a pseudo-natural equivalence,
the assertion for $\Can$ follows from corollary
\ref{cor_fully-faith-2-adjoint}. The same arguments applies
to $\lex.\Can$.
\end{proof}

\begin{remark}
As an application of theorem \ref{th_adj-topos-site} we deduce
the representability of the $2$-limit of any cofiltered system
of topoi. Indeed, let $I$ be a small cofiltered category, and
$T_\bullet:I\to(\sU,\sU')\tdu\Topos$ any pseudo-functor. By
proposition \ref{prop_2-lims-lex-sites}, the pseudo-functor
$\lex.\Can\circ T_\bullet:I\to\lex.\Site$ admits a $2$-limit
$(C,\pi_\bullet)$. By theorem \ref{th_adj-topos-site}(ii)
and proposition \ref{prop_2-adjs-and-lims}(i), the pair
$(\lex.\sT(C),\lex.\sT*\pi_\bullet)$ is a $2$-limit of the
pseudo-functor
$\lex.\sT\circ\lex.\Can\circ T_\bullet:I\to(\sU,\sU')\tdu\Topos$.
But by theorem \ref{th_adj-topos-site}(iii) and corollary
\ref{cor_fully-faith-2-adjoint}, we have a pseudo-natural
equivalence $\omega:\lex.\sT\circ\lex.\Can\isom\one_{\Topos}$,
so $(\lex.\sT(C),(\omega*T_\bullet)\odot(\lex.\sT*\pi_\bullet))$ is
a $2$-limit for $T_\bullet$ (lemma \ref{lem_pseudo-trivial}).
Combining with remark \ref{rem_2-lims-of-sites}, we see that
the $2$-limit of $T_\bullet$ is also a $2$-limit of the induced
pseudo-functor
$$
I\to\sU'\tdu\bCat
\qquad
i\mapsto T_i
\qquad
(\phi:i\to j)\mapsto(T_{\phi*}:T_i\to T_j).
$$
\end{remark}

\subsection{Fibred sites}
\label{sec_morph-of-fibredtop}
The formalism developed in this section and the next one shall
allow us to deal with families of sites or of topoi, indexed by
arbitrary categories. Especially, we will explain how to combine
the members of such a family into a single {\em total site} or
{\em total topos}.

\begin{definition}\label{def_fibred-site}
(i)\ \
Let $\cC$ be a category and $\sV$ a universe. A {\em fibred
site with $\sV$-small fibres over $\cC$} is a datum
$(\cA,p,J_\bullet)$ consisting of a fibration with $\sV$-small
fibre categories
$$
p:\cA\to\cC
$$
and of a topology $J_X$ on the fibre category $\cA_X$,
for every $X\in\Ob(\cC)$, such that the following holds.
For every cleavage $\blambda$ of $p$ and every morphism
$\phi:X\to Y$ in $\cC$, the functor $\sc_\phi:\cA_Y\to\cA_X$
is a morphism of sites $(\cA_X,J_X)\to(\cA_Y,J_Y)$, where
$\sc:\cC^o\to\sV\tdu\bCat$ is the pseudo-functor associated
with $\blambda$.

(ii)\ \
A datum $(\cA,p,J_\bullet)$ as in (i) is a
{\em fibred lex-site with $\sV$-small fibres} if
$(\cA_X,J_X)$ is a lex-site for every $X\in\Ob(\cC)$
and the functors $\sc_\phi$ are left exact for every
morphism $\phi$ of $\cC$ and every cleavage $\blambda$
with associated pseudo-functor $\sc$. As usual, we shall
mostly omit mentioning the universe $\sV$, unless the
omission would be source of ambiguities.

(iii)\ \
Let $(\cA,p,J_\bullet)$ be a fibred site over the category $\cC$.
We endow $\cA$ with the coarsest topology $J_\cA$ such that the
inclusion functor $\cA_X\to\cA$ is continuous for the topologies
$J_X$ and $J_\cA$, for every $X\in\Ob(\cC)$. The {\em total site}
of $(\cA,p,J_\bullet)$ is the resulting site
$$
(\cA,J_\cA).
$$

(iv)\ \
Let $(\cA,p,J_\bullet)$ and $(\cA',p',J'_\bullet)$ be two fibred sites
over $\cC$. A {\em morphism of fibred sites}
$$
\phi:(\cA,p,J_\bullet)\to(\cA',p',J'_\bullet)
$$
is a $\cC$-cartesian functor $\phi:\cA'\to\cA$ whose restriction
$\phi_X:\cA'_X\to\cA_X$ is a morphism of sites
$(\cA_X,J_X)\to(\cA'_X,J'_X)$, for every $X\in\Ob(\cC)$. If
$(\cA,p,J_\bullet)$ and $(\cA',p',J'_\bullet)$ are fibred lex-sites,
we say that $\phi$ is a {\em morphism of fibred lex-sites}, if
it is $\cC$-cartesian, and the restriction
$\phi_X:(\cA_X,J_X)\to(\cA'_X,J'_X)$ is a morphism of lex-sites,
for every $X\in\Ob(\cC)$.
\end{definition}

\begin{remark}\label{rem_fibred-site}
(i)\ \
With the notation of definition \ref{def_fibred-site},
it is clear that if there exists a cleavage of $p$ with
associated pseudo-functor $\sc$, such that $\sc_\phi$ is
a morphism of sites $(\cA_X,J_X)\to(\cA_Y,J_Y)$ for every
morphism $\phi:X\to Y$, then the same holds for every
cleavage of $p$ (and in this case, $p$ is a fibred site).
Likewise, if $\sc_\phi$ is also left exact for every such
$\phi$, then the same holds for every cleavage of $p$
(and then $p$ is a fibred lex-site).

(ii)\ \
Let $(\cA,p,J_\bullet)$ and $(\cA,J_\cA)$ be as in definition
\ref{def_fibred-site}(iii). By lemma \ref{lem_crit-continuity},
we have $\cS\in J_\cA(A)$ for every $A\in\Ob(\cA)$ and every
sieve $\cS\subset\cA/A$ such that
\set\begin{equation}\label{eq-guilde-docs}
\cS\cap(\cA_{pA}/A)\in J_{pA}(A).
\end{equation}

(iii)\ \
In the situation of (ii), suppose moreover that $\cC$ is
small and $(\cA_X,J_X)$ is a $\sU$-site for every $X\in\Ob(\cC)$;
under these assumptions, the fibre category $\cA_X$ has small
$\Hom$-sets for every such $X$, but the same does not necessarily
hold for the category $\cA$. However, for every $A,A'\in\Ob(\cA)$
the set $\Hom_\cA(A,A')$ is essentially small : indeed, let $\blambda$
be a cleavage for $\cA$ and $\sc$ its associated pseudo-functor;
since every morphism $f:A'\to A$ in $\cA$ factors as a morphism
$A'\to\sc_{p(f)}A$ of $\cA_{pA'}$ and the cartesian morphism
$\blambda(A,p(f)):\sc_{p(f)}A\to A$, the assertion follows easily.
Thus, after replacing $\cA$ by an isomorphic category, we may
assume that $\cA$ has small $\Hom$-sets, and in this case we claim
that $(\cA,J_\cA)$ is a $\sU$-site. Indeed, choose for every
$X\in\Ob(\cC)$ a small topologically generating family
$G_X\subset\Ob(p^{-1}X)$ for $(\cA_X,J_X)$; then it is
clear from (ii) that $\bigcup_{X\in\Ob(\cC)}G_X$ is a small
topologically generating family for $(\cA,J_\cA)$.

(iv)\ \
Clearly the fibred sites (resp. the fibred lex-sites) over
$\cC$ with $\sV$-small fibres are the objects of a $2$-category :
$$
\sV\tdu\fib.\Site(\cC)
\qquad
\text{(resp. $\sV\tdu\fib.\lex.\Site$)}
$$
whose $1$-cells are the morphisms of fibred sites (resp. of
fibred lex-sites), and whose $2$-cells $g\Rightarrow g'$ are
the natural $\cC$-transformations $g'\Rightarrow g$ between
the underlying functor of such morphisms, with the obvious
composition laws for $1$-cells and $2$-cells.
We have an obvious forgetful strict pseudo-functor :
$$
\sV\tdu\fib.\Site(\cC)\to{}^o(\sV\tdu\Fib(\cC))^o
\qquad
(\cA,p,J_\bullet)\mapsto(\cA,p).
$$

(v)\ \
Taking into account theorem \ref{th_fundamental-fibrations},
it is easily seen that the pseudo-functor $\cFib_\cC$ induces
strict and strong $2$-equivalences of $2$-categories :
$$
\begin{aligned}
\underline\cFib_\cC&\,:\sPsFun(\cC,\sV\tdu\Site)\to
\sV\tdu\fib.\Site(\cC) \\
\lex.\underline\cFib_\cC&\,:\sPsFun(\cC,\sV\tdu\lex.\Site)\to
\sV\tdu\fib.\lex.\Site(\cC).
\end{aligned}
$$
Namely, we have an obvious forgetful strict pseudo-functor
$$
\Phi:\sV\tdu\Site\to{}^o(\sV\tdu\bCat)^o
\qquad
(\cA,J)\mapsto\cA
$$
whence a pseudo-functor $\sPsFun(\cC,\Phi):
\sPsFun(\cC,\sV\tdu\Site)\to\sPsFun(\cC,{}^o(\sV\tdu\bCat)^o)$
(see remark \ref{rem_remains-of-the-day}(i)). Then, since
${}^o\cC^o=\cC^o$, we have a strict isomorphism of $2$-categories :
$$
\sPsFun(\cC,{}^o(\sV\tdu\bCat)^o)\isom
{}^o\sPsFun(\cC^o,\sV\tdu\bCat)^o
$$
(see \eqref{subsec_opposite-mods}) and a strict and strong
$2$-equivalence
$$
{}^o\cFib_\cC^o:{}^o\sPsFun(\cC^o,\sV\tdu\bCat)^o
\isom{}^o(\sV\tdu\Fib(\cC))^o.
$$
The composition
$\sPsFun(\cC,\sV\tdu\Site)\to{}^o(\sV\tdu\Fib(\cC))^o$
of these pseudo-functors assigns to every pseudo-functor
$$
\underline\cA:\cC\to\sV\tdu\Site
\qquad
X\mapsto(\cA_X,J_X)
$$
the fibration $\cA':=\cFib({}^o(\Phi\circ\underline\cA)^o)$ whose
fibre category $\cA'_X$ over every $X\in\Ob(\cC)$ is naturally
identified with $\cA_X$. Then $\underline\cFib(\underline\cA)$
is the fibred site $(\cA',J'_\bullet)$ where $J'_X$ is the topology on
$\cA'_X$ corresponding to $J_X$ under this identification
$\cA'_X\isom\cA_X$, for every $X\in\Ob(\cC)$.

(vi)\ \
By remark \ref{rem_strict-strong}, the pseudo-functor
$\underline\cFib_\cC$ admits a strict and strong pseudo-inverse
$$
\underline\sc^\bullet:\sV\tdu\fib.\Site\to
\sPsFun(\cC,\sV\tdu\Site)
\qquad
\underline\cA:=(\cA,p,J_\bullet)\mapsto\underline\sc^{\underline\cA}
$$
which, to every fibred site $(\cA,p,J_\bullet)$ over $\cC$,
assigns the pseudo-functor $\sc^\cA:\cC^o\to\sV\tdu\bCat$
associated with a cleavage of $\cA$, and for every $X\in\Ob(\cC)$,
endows the category $\sc^\cA_X=\cA_X$ with the topology $J_X$.
Obviously $\underline\sc^\bullet$ restricts to a strict and
strong pseudo-inverse for $\lex.\underline\cFib_\cC$.
\end{remark}

\sset\subsubsection{}\label{subsec_sheaves-on-tot-site}
Let $\cC$ be a category, $(\cA,p,J_\bullet)$ a fibred lex-site
over $\cC$, and $(\cA,J_\cA)$ its total site. For every
$X\in\Ob(\cC)$ denote by $i_X:\cA_X\to\cA$ the inclusion
functor. We have :

\begin{proposition}\label{prop_sheaves-on-tot-site}
In the situation of \eqref{subsec_sheaves-on-tot-site}, the
following holds :
\begin{enumerate}
\item
The functor $i_X$ commutes with fibre products and equalizers
for every $X\in\Ob(\cC)$.
\item
For every $A\in\Ob(\cA)$ and every sieve $\cS\subset\cA/A$,
we have $\cS\in J_\cA(A)$ if and only if $\cS$ satisfies condition
\eqref{eq-guilde-docs}.
\item
The functor $i_X$ is both continuous and cocontinuous for
the topologies $J_X$ and $J_\cA$, for every $X\in\Ob(\cC)$.
\item
A presheaf $F$ on $\cA$ is a sheaf for the topology $J_\cA$ if
and only if $i^\wedge_XF$ is a sheaf on $(\cA_X,J_X)$ for every
$X\in\Ob(\cC)$.
\item
Suppose moreover that $\cC$ is small, and $(\cA_X,J_X)$ is a
$\sU$-site for every $X\in\Ob(\cC)$. Then we have an essentially
commutative diagram of categories :
$$
{\spreaddiagramcolumns{+20pt}\diagram
\cA^\wedge \ar[r]^-{i^\wedge_X} \ar[d]_{(-)^a} &
\cA_X^\wedge \ar[d]^{(-)^a} \\
(\cA,J_\cA)^\sim \ar[r]^-{i^\sim_{X*}} & (\cA_X,J_X)^\sim
\enddiagram}
\qquad
\text{for every $X\in\Ob(\cC)$}.
$$
\end{enumerate}
\end{proposition}
\begin{proof} Choose a cleavage $\blambda$ for $p$ and let
$\sc:\cC^o\to\sV\tdu\bCat$ be the associated pseudo-functor
(for some universe $\sV$ such that $\cA_X$ is a $\sV$-small
category for every $X\in\Ob(\cC)$).

(i): Let $A'\xrightarrow{f'}A\xleftarrow{f''}A''$ be two
morphisms in $\cA_X$; pick $P\in\Ob(\cA_X)$ representing
$A'\times_AA''$ in $\cA_X$, and let
$A'\xleftarrow{q'}P\xrightarrow{q''}A''$ be a universal cone.
Let then $B\in\Ob(\cA)$ be any object and
$A'\xleftarrow{g'}B\xrightarrow{g''}A''$ two morphisms in $\cA$
such that $f'\circ g'=f''\circ g''$. It follows easily that
$\phi:=p(g')=p(g'')$, so $g'$ (resp. $g''$) factors through
a morphism $h':B\to\sc_\phi A'$ (resp. $h'':B\to\sc_\phi A''$)
in $\cA_{pB}$ and the cartesian morphism
$\blambda(A',\phi):\sc_\phi A'\to A'$
(resp. $\blambda(A'',\phi):\sc_\phi A''\to A''$), and since
$\sc_\phi$ is left exact, we get a unique morphism
$h:B\to\sc_\phi P$ in $\cA_{pB}$ such that
$\sc_\phi(q')\circ h=h'$ and $\sc_\phi(q'')\circ h=h''$. Then
$k:=\blambda(\phi,P)\circ h:B\to P$ is the unique morphism
in $\cA$ such that $q'\circ k=g'$ and $q''\circ k=g''$. This
shows that $P$ represents the fibre product of $A'$ and $A''$
over $A$ in $\cA$. Similarly one shows the assertion for
equalizers : the details shall be left to the reader.

(ii): In view of (i) and lemma \ref{lem_crit-continuity}, it
suffices to check that the family $\cJ$ of sieves verifying
condition \eqref{eq-guilde-docs} yields a topology on $\cA$.
Now, condition (c) of definition \ref{def_topology}(i) obviously
holds for $\cJ$. Next, let $\cS\subset\cA/A$ be a sieve
fulfilling condition \eqref{eq-guilde-docs}, and $f:A'\to A$
a morphism of $\cA$; we set $\phi:=p(f)$ and factor $f$ through
a morphism $f':A'\to A'':=\sc_\phi A$ in $p^{-1}A'$ and the
cartesian morphism $f'':=\blambda(A,\phi):A''\to A$. Then
let $\cS':=\cS\times_Af''$, and notice that if
$(g_\lambda:B_\lambda\to A~|~\lambda\in\Lambda)$ generates
$\cS_{|A}:=\cS\cap(\cA_{pA}/A)$, then the family
$(\sc_\phi(g_\lambda):\sc_\phi B_\lambda\to A''~|~\lambda\in\Lambda)$
lies in $\cS'$, and the latter family covers $A''$ for the
topology $J_{pA''}$, since $\sc_\phi$ is continuous for the
topologies $J_{pA}$ and $J_{pA''}$, and since $\cS_{|A}$ covers $A$
for the topology $J_{pA}$. Then $\cS\times_Af=\cS'\times_{A''}f'$
contains the family
$(A'\times_{A''}\sc_\phi B_\lambda\to A'~|~\lambda\in\Lambda)$,
which covers $A'$ for the topology $J_{pA'}$. This shows that
condition (a) of definition \ref{def_topology}(i) holds for
$\cJ$. Lastly, let $\cS\subset\cA/A$ be an element of $\cJ$,
and consider another sieve $\cT\subset\cA/A$ such that for
every $(f:A'\to A)\in\Ob(\cS)$ the sieve $\cT\times_Af$ lies
in $\cJ$. Let $\cS_{|A}\subset\cA_{pA}/A$ be as in the
foregoing, and define likewise $\cT_{|A}$; then
$\cS_{|A}\in J_{pA}(A)$ and it is easily seen that
$$
\cT_{|A}\times_Af=(\cT\times_Af)\cap(\cA_{pA'}/A')
\qquad
\text{for every $(f:A'\to A)\in\Ob(\cS_{|A})$}.
$$
Hence $\cT_{|A}\in J_{pA}(A)$, so $\cT\in\cJ$;
this shows that $\cJ$ fulfills condition (b) of
definition \ref{def_topology}(i).

(iii) follows immediately from (ii).

(v): Pick a universe $\sV$ containing $\sU$, and such
that $\cA$ is $\sV$-small; in view of remarks \ref{rem_U-site}(i)
and \ref{rem_fibred-site}(iii) we are easily reduced to checking
the essential commutativity of the diagram :
$$
{\spreaddiagramrows{-5pt}\spreaddiagramcolumns{+20pt}\diagram
\cA^\wedge_\sV \ar[r]^-{i^\wedge_{X,\sV}} \ar[d]_{(-)^+} &
\cA_{X,\sV}^\wedge \ar[d]^{(-)^+} \\
\cA^\wedge_\sV \ar[r]^-{i^\wedge_X} & \cA^\wedge_{X,\sV}
\enddiagram}
\qquad
\text{for every $X\in\Ob(\cC)$}
$$
where the functors $(-)^+$ are defined as in
\eqref{subsec_asso-topoi}. We may then replace $\sU$ by $\sV$
and assume that $\cA$ is a small category. For every
$X\in\Ob(\cC)$ and every $A\in\Ob(\cA_X)$ we endow $J_X(A)$
and $J_\cA(A)$ with the ordering induced by inclusion of sieves,
and consider the map
$$
J_X(A)\to J_\cA(A)
\qquad
\cS\mapsto\cS^*
$$
that assigns to every $\cS\in J_X(A)$ the sieve $\cS^*\in J_\cA(A)$
generated by $\Ob(\cS)\subset\Ob(\cA/A)$. This map is obviously
injective, and according to (ii), its image is a cofinal subset
in the opposite ordered set $J_\cA(A)^o$. Moreover, for every
morphism $f:A'\to A$ in $\cA_X$, clearly we have
$$
(\cS\times_Af)^*=\cS^*\times_Af
\qquad
\text{for every $\cS\in J_X(A)$}.
$$
Furthermore, for every such $\cS$ and every presheaf $F$ on
$\cA_X$, we have a natural identification of
$F(\cS):=\Hom_{\cA_X^\wedge}(h_\cS,F)$ with the equalizer of the
two natural maps (see \eqref{subsec_interpret-descent})
\set\begin{equation}\label{eq_describe-F(S)}
\Equal\Bigl( \prod_{(B\to A)\in\Ob(\cS)}
\xymatrix{\!\!\!\!\!F(B)\ar@<-.5ex>[r] \ar@<.5ex>[r] &}
\!\!\!\!\!\!\!\!\!\!\!\!
\prod_{(B\to A\leftarrow B')\in\Ob(\cS)\times\Ob(\cS)}\!\!\!\!\!\!\!\!F(B\times_AB')
\Bigr)
\end{equation}
and for every inclusion $\cS'\subset\cS$ of sieves covering $A$,
the map $F(\cS)\to F(\cS')$ corresponds, under this identification,
to the restriction of the projection
\set\begin{equation}\label{eq_restrict-this}
\prod_{(B\to A)\in\Ob(\cS)}\!\!\!\!\!F(B)\to\!\!\!\!\!
\prod_{(B\to A)\in\Ob(\cS')}\!\!\!\!\!F(B).
\end{equation}
Recall that
$$
F^+(A):=\colim_{\cS\in J_X(A)^o}F(\cS)
$$
and for every morphism $f:A'\to A$ in $\cA_X$, the map
$F^+(f):F^+(A)\to F^+(A')$ is the colimit of the maps
$F(\cS)\to F(\cS\times_Af)\to F^+(A')$ induced by the
projection $\cS\times_Af\to\cS$, for every $\cS\in J_X(A)^o$.
Now, if $F=i^\wedge_XG$ for a presheaf $G$ on $\cA$, it follows
from \eqref{subsec_interpret-descent} that the equalizer
\eqref{eq_describe-F(S)} is also naturally identified with
$G(\cS^*)$, and likewise, for every inclusion of sieves
$\cS'\subset\cS$ in $J_X(A)$, the map $G(\cS^*)\to G(\cS'^*)$
corresponds to the restriction of \eqref{eq_restrict-this},
under this same identification. After taking colimits, we thus
get an isomorphism
$$
(i^\wedge_XG)^+(A)\isom G^+(A)
$$
natural in both the presheaf $G$ and the object
$A\in\Ob(\cA_X)$. The assertion follows.

(iv): The condition is obviously necessary. For the converse,
let $F$ be a presheaf on $\cA$ such that $i^\wedge_X(F)$ is a
sheaf on $(\cA_X,J_X)$ for every $X\in\Ob(\cC)$; denote by
$F^a$ the sheaf on $(\cA,J_\cA)$ associated with $F$, and
by $j:F\to F^a$ the natural map of presheaves. For every
$X\in\Ob(\cC)$ we get a commutative diagram of presheaves on
$\cA_X$ :
$$
\xymatrix@R-5pt{ i^\wedge_X(F) \ar[rr]^-{i^\wedge_X(j)} \ar[d] & &
i^\wedge_X(F^a) \ar[d] \\
(i^\wedge_XF)^a \ar[rr]^-{(i^\wedge_X(j))^a} & & (i^\wedge_XF^a)^a
}$$
whose vertical arrows are isomorphisms, by assumption. On the
other hand, in view of (v), the morphism $(i^\wedge_X(j))^a$ is
naturally identified with
$i^\wedge_X(j^a):i^\wedge_X(F^a)\to i^\wedge_X((F^a)^a)$; now, $j^a$
is an isomorphism, hence the same holds for $i^\wedge_X(j^a)$,
and then also for $i^\wedge_X(j)$, for every $X\in\Ob(\cC)$.
But this means that $j$ is already an isomorphism, {\em i.e.}
$F$ is a sheaf on $(\cA,J_\cA)$.
\end{proof}

\begin{example}\label{ex_base-with-fin-obj}
In the situation of \eqref{subsec_sheaves-on-tot-site}, suppose
that $\cC$ admits a final object $X_0$. Then the functor $i_{X_0}$
is a morphism of sites $(\cA,J_\cA)\to(\cA_{X_0},J_X)$. Indeed,
let $A_0$ be any final object of $\cA_{X_0}$; it is easily
seen that $A_0$ is also a final object of $\cA$, and taking
into account propositions \ref{prop_commute-criteria}(i) and
\ref{prop_sheaves-on-tot-site}(i) we deduce that $i_{X_0}$ is
left exact. Then the assertion follows from proposition
\ref{prop_sheaves-on-tot-site}(iii) and example \ref{ex_simple-case}.
\end{example}

\begin{proposition}\label{prop_actually-morph-of-sites}
Let $u:(\cA',p',J'_\bullet)\to(\cA,p,J_\bullet)$ be a morphism
of fibred lex-sites over a category $\cC$, and denote by
$(\cA,J_\cA)$, $(\cA',J_{\cA'})$ the corresponding total sites.
Then $u$ is a morphism of sites $(\cA',J_{\cA'})\to(\cA,J_\cA)$.
\end{proposition}
\begin{proof} Proposition \ref{prop_sheaves-on-tot-site}(i,ii)
easily implies that $u$ is continuous for the topologies $J_\cA$
and $J_{\cA'}$, so it remains only to check that the fibration
$\ss:\cA'/u\cA\to\cA'$ is locally cofiltered for the topology
$J_{\cA'}$ (proposition \ref{prop_morph-of-sites}).
Thus, let $X\in\Ob(\cC)$ and $A'\in\Ob(\cA'_X)$; the fibration
$\ss_{(X)}:\cA'_X/u_X\cA_X\to\cA'_X$ is locally cofiltered, by
proposition \ref{prop_morph-of-sites}, hence there exists
a covering family
$g_\bullet:=(g_\lambda:A'_\lambda\to A'~|~\lambda\in\Lambda)$
relative to the topology $J'_X$ such that
$\Ob(\ss_{(X)}^{-1}A'_\lambda)\neq\emptyset$ for every
$\lambda\in\Lambda$; but the family $g_\bullet$ covers $A'$
also relative to the topology $J_{\cA'}$ (proposition
\ref{prop_sheaves-on-tot-site}(ii)), so $\ss$ fulfills
condition (a) of definition \ref{def_locally-cofiltered}.
To check condition (b), consider morphisms $f_i:A'\to uA_i$
in $\cA'$ for $i=1,2$, and let $\phi_i:X\to Y_i$ be the
image of $f_i$ in $\cC$ for $i=1,2$. Pick any cleavage
for $\cA$, and let $\sc$ be its associated pseudo-functor;
then $f_i=u(f''_i)\circ f'_i$ where $f'_i:A'\to u(\sc_{\phi_i}A_i)$
is a morphism in $\cA'_X$ and $f''_i:\sc_{\phi_i}A_i\to A_i$
is a cartesian morphism of $\cA$, for $i=1,2$.
Since $\ss_{(X)}$ is locally cofiltered, we may then find
a covering family $g_\bullet$ as in the foregoing, and for
every $\lambda\in\Lambda$ a morphism
$f_\lambda:A'_\lambda\to uA_\lambda$ in $\cA'_X$ and morphisms
$h_{i,\lambda}:A_\lambda\to\sc_{\phi_1}A_1$ in $\cA_X$ for $i=1,2$,
such that $u(h_{i,\lambda})\circ f_\lambda=f'_i\circ g_\lambda$.
Then we get morphisms in $\cA/u\cA'$
$$
(A'_\lambda/f''_i\circ h_{i,\lambda}):(A'_\lambda,f_\lambda)\to
(A'_\lambda,f_i\circ g_\lambda)
\qquad
\text{for $i=1,2$ and every $\lambda\in\Lambda$}
$$
whence the contention. Lastly, let $(f_i:A'\to uA_i~|~i=1,2)$
be a pair of morphisms in $\cA'$, and $h_1,h_2:A_1\to A_2$
two morphisms of $\cA$ with $u(h_i)\circ f_1=f_2$ for $i=1,2$.
Let $\phi_i:X\to Y_i$ be the image of $f_i$ in $\cC$ for $i=1,2$;
there exist morphisms $h'_1,h'_2:\sc_{\phi_1}A_1\to\sc_{\phi_2}A_2$
in $\cA_X$ and cartesian morphisms $f'_i:\sc_{\phi_i}A_i\to A_i$
in $\cA$ such that $h_i\circ f'_1=f'_2\circ h'_i$ for $i=1,2$.
Since $\ss_{(X)}$ is locally cofiltered, we may then find
$g_\bullet$ as in the foregoing, and for every $\lambda\in\Lambda$
a morphism $f_\lambda:A'_\lambda\to uA_\lambda$ in $\cA'_X$ and
a morphism $h_\lambda:A_\lambda\to\sc_{\phi_1}A_1$ in $A_X$ such
that $h'_1\circ h_\lambda=h'_2\circ h_\lambda$. It follows that
$h_1\circ f'_1\circ h_\lambda=h_2\circ f'_1\circ h_\lambda$,
which shows that condition (c) holds as well for $\ss$.
\end{proof}

\begin{remark}\label{rem_actually-morph-of-sites}
Let $\sV$ be any universe. Proposition
\ref{prop_actually-morph-of-sites} says that the rule that
assigns to every fibred lex-site over $\cC$ with $\sV$-small
fibres its total site defines a strict pseudo-functor
$$
\totSite:\sV\tdu\fib.\lex.\Site(\cC)\to\sV\tdu\Site
$$
that assigns to every morphism of lex-sites the underlying
functor, and is the identity on $2$-cells.
\end{remark}

\sset\subsubsection{}\label{subsec_pullback-fib-site}
Let $(\cA,p,J_\bullet)$ be a fibred site over
a category $\cC$, and $u:\cC'\to\cC$ any functor. Set
$\cA':=\Fib(u)^*\cA$, and recall that for every $X\in\Ob(\cC')$
the natural projection $\pi:\cA'\to\cA$ restricts to an isomorphism
of categories $\pi_X:\cA'_X\isom\cA_{uX}$. Then we get a fibred site
denoted
$$
\cC'\times_{(\cC,u)}(\cA,p,J_\bullet)
\qquad\text{or more simply}\qquad
\cC'\times_\cC(\cA,p,J_\bullet)
$$
whose underlying fibration is $p':\cA'\to\cC'$ and where the
topology $J_X$ on $\cA'_X$ is induced by $J_X$ via the
isomorphism $\pi_X$, for every $X\in\Ob(\cC')$. Obviously,
if $(\cA,p,J_\bullet)$ is a fibred lex-site, the same holds
for $\cC'\times_\cC(\cA,p,J_\bullet)$, and in this case,
proposition \ref{prop_sheaves-on-tot-site}(ii,iv) easily
implies that $\pi$ is a continuous and cocontinuous functor
for the respective total sites $(\cA,J_\cA),(\cA',J_{\cA'})$.

\sset\subsubsection{}\label{subsec_morph-fibsites-from-nat-tr}
In the situation of \eqref{subsec_pullback-fib-site}, let
$v:\cC'\to\cC$ be another functor, and $\beta:u\Rightarrow v$
any natural transformation. After choosing a cleavage $\blambda$
for the fibration $\cA$, we attach to $\beta$ a $\cC'$-cartesian
functor $\Fib(\beta)^*_\cA:\Fib(v)^*\cA\to\Fib(u)^*\cA$, as in
\eqref{subsec_complete-Fib}. Explicitly, let $\sc$ be the
pseudo-functor associated with $\blambda$; then for every
$X\in\Ob(\cC')$ the restriction $\cA_{vX}\to\cA_{uX}$ of
$\Fib(\beta)^*_\cA$ is given by the functor $\sc_{\beta_X}$. Hence,
$\Fib(\beta)^*_\cA$ is a morphism of fibred sites denoted :
$$
\beta\times_\cC(\cA,p,J_\bullet):
\cC'\times_{(\cC,u)}(\cA,p,J_\bullet)\to
\cC'\times_{(\cC,v)}(\cA,p,J_\bullet).
$$
It is easily seen that $\beta\times_\cC(\cA,p,J_\bullet)$ is
independent, up to isomorphisms of morphisms of fibred sites,
of the choice of cleavage $\blambda$ : the details shall be
left to the reader.

\sset\subsubsection{}\label{subsec_box-a-louer}
Let $\cC,\cD$ be two small categories, $t:\cD\to\cC$ a functor,
and $u:(\cA_0,p_0,J_{0,\bullet})\to(\cA_1,p_1,J_{1,\bullet})$ a
morphism of fibred lex sites over $\cC$. Set
$$
(\cB_i,q_i,J'_{i,\bullet}):=\cD\times_\cC(\cA_i,p_i,J_{i,\bullet})
\qquad
\text{for $i=0,1$}
$$
(notation of \eqref{subsec_pullback-fib-site}). Clearly $u$ induces
a morphism of fibred lex sites
$$
v:=\cD\times_\cC u:(\cB_0,q_0,J'_{0,\bullet})\to(\cB_1,q_1,J'_{1,\bullet}).
$$
Then $u$ and $v$ are also morphisms of the respective total sites
$u:(\cA_0,J_0)\to(\cA_1,J_1)$ and $v:(\cB_0,J'_0)\to(\cB_1,J'_1)$
(proposition \ref{prop_actually-morph-of-sites}). For $i=0,1$,
let $\pi_i:\cB_i\to\cA_i$ be the projection. We deduce a commutative
diagram of categories :
$$
\cE \qquad : \qquad
{\diagram \cA_0^\sim \ar[r]^-{\tilde u_*} \ar[d]_{\tilde\pi_{0*}} &
\cA_1^\sim \ar[d]^{\tilde\pi_{1*}} \\
\cB_0^\sim \ar[r]^-{\tilde v_*} & \cB_1^\sim
\enddiagram}
\qquad\qquad
$$
and notice that each functor of the diagram $\cE$ admits a left
adjoint; after choosing an adjunction for each of the resulting
adjoint pairs of functors, we can view $\cE$ as a square of
links, oriented by the identity $\one_{\tilde\pi_{1*}\circ\tilde u_*}$
(see \eqref{subsec_base-change-map}).

\begin{proposition}\label{prop_Devos}
In the situation of \eqref{subsec_box-a-louer}, suppose that
$(\cA_{i,X},J_{i,X})$ is a $\sU$-site for every $X\in\Ob(\cC)$
and $i=0,1$. Then the base change transformation
$$
\Upsilon(\cE):
\tilde v{}^*\circ\tilde\pi_{1*}\to\tilde\pi_{0*}\circ\tilde u{}^*
$$
is an isomorphism of functors.
\end{proposition}
\begin{proof} We consider first the oriented diagram of
categories :
$$
\cE' \qquad :\qquad
{\diagram
\cA_0^\wedge \ar[rrr]^-{u^\wedge} \ar[ddd]_{\pi_0^\wedge}
& & \dltwocell\omit{^\ } & \cA^\wedge_1 \ar[ddd]^{\pi_1^\wedge} \\
& \cA_0^\sim \ar[r]^-{\tilde u_*} \dltwocell\omit{^\ }
\drtwocell\omit{} \ar[d]_{\tilde\pi_{0*}} \ar[ul]_-{i_{\cA_0}} &
\cA_1^\sim \ar[d]^{\tilde\pi_{1*}} \ar[ur]^-{i_{\cA_1}}
\drtwocell\omit{_\ } \\
& \cB_0^\sim \ar[r]^-{\tilde v_*} \ar[dl]_-{i_{\cB_0}}
\drtwocell\omit{^\ } & \cB_1^\sim \ar[dr]^-{i_{\cB_1}} & \\
\cB_0^\wedge \ar[rrr]^-{v^\wedge} & & & \cB_1^\wedge
\enddiagram}
\qquad\qquad$$
whose diagonal arrows are the inclusion functors, and all
whose orientations are identities. Moreover, each functor
in $\cE'$ still admits a left adjoint, so as usual $\cE'$
can be regarded as an oriented diagram of links. By adding
a further identity orientation for the ``front face'' we
obtain a cubical diagram as in
\eqref{subsec_transfer-base-change}, which obviously
commutes on $2$-cells, except that the orientations do
not agree with those of the corresponding diagram in
\eqref{subsec_transfer-base-change}, nor with those of its
variant considered in remark \ref{rem_transit-base-change}(iii).
However, the induced diagram ${}^\dagger\cE'$ of adjoint squares
as in \eqref{subsec_base-change-map} still commutes on $2$-cells,
due to proposition \ref{prop_isoms-of-links}, and moreover its
orientations agree with those of the variant from remark
\ref{rem_transit-base-change}(iii). More precisely, the
left and right faces of ${}^\dagger\cE'$ are oriented by
${}^o(\one_{i_{\cB_1}\circ\tilde v_*})^\dagger$ and respectively
${}^o(\one_{i_{\cA_1}\circ\tilde u_*})^\dagger$, and the top and
bottom faces are oriented by
${}^o(\one_{i_{\cB_1}\circ\tilde\pi_{1*}})^\dagger$ and respectively
${}^o(\one_{i_{\cB_0}\circ\tilde\pi_{0*}})^\dagger$.

\begin{claim}\label{cl_birthday-tomorrow}
The base change transformations
$$
\Upsilon(\one_{i_{\cB_1}\circ\tilde\pi_{1*}}):
(-)^a_{\cB_1}\circ\pi^\wedge_1\to\tilde\pi_{1*}\circ(-)^a_{\cA_1}
\qquad\text{and}\qquad
\Upsilon(\one_{i_{\cB_0}\circ\tilde\pi_{0*}}):
(-)^a_{\cB_0}\circ\pi^\wedge_0\to\tilde\pi_{0*}\circ(-)^a_{\cA_0}
$$
are isomorphisms of functors.
\end{claim}
\begin{pfclaim} We check the assertion for
$\Upsilon(\one_{i_{\cB_0}\circ\tilde\pi_{0*}})$; the same argument
shall apply also to $\Upsilon(\one_{i_{\cB_1}\circ\tilde\pi_{1*}})$.
Since $i_{\cB_0}$ is a fully faithful functor, the counit
$\eps^{\cB_0}$ of the chosen adjunction for the pair
$((-)^a_{\cB_0},i_{\cB_0})$ is an isomorphism (proposition
\ref{prop_fullfaith-adjts}(iii)), so we are reduced to
checking that $((-)^a_{\cB_0}\circ\pi^\wedge_0)*\eta^{\cA_0}$
is an isomorphism, where $\eta^{\cA_0}$ and $\eps^{\cA_0}$
are the unit and counit of an adjunction for the pair
$((-)^a_{\cA_0},i_{\cA_0})$. But recall that there exists
an isomorphism of functors :
$(-)^a_{\cB_0}\circ\pi_0^\wedge\isom\tilde\pi_{0*}\circ(-)^a_{\cA_0}$
(corollary \ref{cor_two-U-sites}(iii)), so it suffices
to show that $(-)^a_{\cA_0}*\eta^{\cA_0}$ is an isomorphism.
However, $\eps^{\cA_0}$ is an isomorphism, since $i_{\cA_0}$
is fully faithful, hence the contention follows from the
triangular identities of \eqref{subsec_adj-pair}.
\end{pfclaim}

In light of claim \ref{cl_birthday-tomorrow}, proposition
\ref{prop_opp-links-and-base-ch}, and remark
\ref{rem_transit-base-change}(i,iii), we are then reduced
to showing that the base change transformation
$$
\Upsilon(\one_{\pi_1^\wedge\circ u^\wedge}):
v_!\circ\pi_1^\wedge\to\pi_0^\wedge\circ u_!
$$
is an isomorphism of functors. To this aim, recall that for
every $A_0\in\Ob(\cA_0)$ and every presheaf $F$ on $\cA_1$,
the set $u_!F(A_0)$ represents the colimit of the functor
$$
F\circ\st^o_{A_0}:(A_0/u\cA_1)^o\to\Set
$$
where $\st_{A_0}$ denotes the usual target functor as in
\eqref{subsec_fibreovercat}. Thus, $u_!F(A_0)$ is the
set of equivalence classes $[s,\phi]_u$ of pairs $(s,\phi)$,
where $\phi:A_0\to uA_1$ is a morphism of $\cA_0$, and
$s\in FA_1$. For every morphism $\psi:A_0\to A'_0$ in
$\cA_0$, the map $u_!F(\psi)$ is given by the rule :
$[s,\phi]_u\mapsto[s,\phi\circ\psi]_u$ for every such
$[s,\phi]_u$. Moreover, for every morphism $f:F\to F'$
of presheaves on $\cA_0$, the morphism $u_!f:u_!F\to u_!F'$
is given by the rule : $[s,\phi]_u\mapsto[f_{A_1}(s),\phi]_u$
for every such $[s,\phi]_u$. We have a standard adjunction
for the pair $(u_!,u^\wedge)$ whose unit and counit
$(\eta^u,\eps^u)$ are given by the following rules,
for every presheaf $F$ on $\cA_1$ and $F'$ on $\cA_0$ :
$$
\begin{aligned}
&\eta^u_{F,A_1}(s):=[s,\one_{uA_1}]_u\in u^\wedge u_!F(A_1)
&\quad&\text{for every $A_1\in\Ob(\cA_1)$ and every $s\in FA_1$} \\
&\eps^u_{F',A_0}[s,A_0\xrightarrow{\phi}uA'_1]_u:=F'\phi(s)\in F'A_0
&\quad&\text{for every $s\in F(uA'_1)$}.
\end{aligned}
$$
A similar description applies to the functor $v_!$
and the unit and counit of its standard adjunction.
Thus, for every such $F$, the morphism
$(v_!\pi_1^\wedge*\eta^u)_F:
v_!\pi^\wedge_1F\to v_!\pi_1^\wedge u^\wedge u_!F$ is given by the
rule:
$$
[s,B_0\xrightarrow{\phi}vB_1]_v\mapsto[[s,\one_{u\pi_1B_1}]_u,\phi]_v
\qquad
\text{for every $s\in F\pi_1B_1$}
$$
and the morphism $(\eps^v*\pi_0^\wedge u_!)_F:
v_!v^\wedge\pi_0^\wedge u_!F\to\pi_0^\wedge u_!F$ is given by the rule :
$$
[[t,\pi_0 vA'_1\xrightarrow{\psi}uA''_1]_u,A_0\xrightarrow{\phi}vA'_1]_v
\mapsto
\pi^\wedge_0u_!F(\phi)[t,\pi_0 vA'_1\xrightarrow{\psi}uA''_1]_u=
[t,\psi\circ\pi_0(\phi)]_u
$$
for every $t\in FA''_1$. Summing up,
$\Upsilon(\one_{\pi_1^\wedge\circ u^\wedge})$ is then given by the rule :
$$
[s,B_0\xrightarrow{\phi}vB_1]_v\mapsto[s,\pi_0(\phi)]_u
\qquad
\text{for every $s\in F\pi_1B_1$}.
$$
Consider then the functor
$$
\rho^{B_0}:(B_0/v\cB_1)^o\to(\pi_0B_0/u\cA_1)^o
\qquad
(B_0\xrightarrow{\phi}vB_1)\mapsto
(\pi_0B_0\xrightarrow{\pi_0(\phi)}\pi_0vB_1=u\pi_1B_1)
$$
and notice that $\st^o_{\pi_0B_0}\circ\rho^{B_0}=\pi^o_1\circ\st^o_{B_0}$,
so that $\rho^{B_0}$ induces a morphism of presheaves on
the category $\Set^o$, as in (the dual of) remark
\ref{rem_wishful}(i) :
$$
\colim_{(\pi_0B_0/u\cA_1)^o}F\circ\st^o_{\pi_0B_0}\to
\colim_{(\cB_0/v\cB_1)^o}F\circ\pi^o_1\circ\st^o_{B_0}
$$
which corresponds to a map
$\omega:v_!\pi^\wedge_1F(B_0)\to\pi_0^\wedge u_!F(B_0)$, and
by a direct inspection it is easily seen that
$\omega=\Upsilon(\one_{\pi_1^\wedge\circ u^\wedge})_{F,B_0}$.
In view of proposition \ref{prop_MacL-cofinal}, we are then
reduced to showing :

\begin{claim} For every $X\in\Ob(\cB_0)$ we have :

(i)\ \
The categories $X/v\cB_1$ and $\pi_0X/u\cA_1$ are cofiltered.

(ii)\ \
The functor $\rho^X$ is cofinal.
\end{claim}
\begin{pfclaim}[](i): We show that $X/v\cB_1$ is cofiltered;
the same argument shall apply to $\pi_0X/u\cA_1$ as well.
Indeed, let us check that $X/v\cB_1$ is not the empty
category : to this aim, let $v_{q_0X}:\cB_{1,q_0X}\to\cB_{0,q_0X}$
be the restriction of $v$ to the fibre categories over $q_0X$,
so it suffices to check that $X/v_{q_0X}\cB_{1,q_0X}$ is not empty;
but $\cB_{1,q_0X}$ is finitely complete, hence it admits a final
object $Y_0$, and $vY_0$ is a final object of $\cB_{0,q_0X}$,
since $v$ is left exact. The unique morphism $X\to vY_0$
is then an object of $X/v_{q_0X}\cB_{1,q_0X}$.

Next, consider two objects $(Y_1,\phi_1),(Y_2,\phi_2)$ of
$X/v\cB_1$; for $i=1,2$ we may find cartesian morphisms
$\beta_i:Y'_i\to Y_i$ in $\cB_1$ such that
$q_1(\beta_i)=q_0(\phi_i)$, and then there exists a unique
morphism $\phi'_i:X\to vY'_i$ in the fibre category
$\cB_{0,q_0X}$ such that $v(\beta_i)\circ\phi'_i=\phi_i$. The
pair $(Y'_i,\phi'_i)$ is an object of $X/v\cB_1$, and $\beta_i$
is a morphism $(Y'_i,\phi'_i)\to(Y_i,\phi_i)$ for $i=1,2$.
Since $\cB_{1,q_0X}$ is finitely complete, the product
$Z:=Y'_1\times Y'_2$ is representable, and since $v_{q_0X}$ is
left exact, $vZ$ represents $vY'_1\times vY'_2$, whence a
morphism $\psi':X\to vZ$ in $\cB_{0,q_0X}$ whose composition
with the projection $v(p_i):vZ\to vY'_i$ agrees with $\phi'_i$,
for $i=1,2$. Thus, $(Z,\psi')\in\Ob(X/v\cB_1)$ and we have
morphisms $(Z,\psi')\xrightarrow{X/p_i}(Y'_i,\phi'_i)
\xrightarrow{X/\beta_i}(Y_i,\phi_i)$ for $i=1,2$.

Lastly, let $X/\alpha_1,X/\alpha_2:(Y_1,\phi_1)\to(Y_2,\phi_2)$
be two morphisms of $X/v\cB_1$. We find cartesian morphisms
$\beta_i:Y'_i\to Y_i$ in $\cB_1$ such that $q_1(\beta_i)=q_0(\phi_i)$
for $i=1,2$, and notice that
$q_1(\alpha_1\circ\beta_1)=q_0(\phi_2)=q_1(\alpha_2\circ\beta_1)$.
Then there exist unique morphisms $\alpha'_i:Y'_1\to Y'_2$ in
$\cB_{1,q_0X}$ such that $\beta_2\circ\alpha'_i=\alpha_i\circ\beta_1$
for $i=1,2$. Since $\cB_{1,q_0X}$ is finitely complete, the equalizer
of the pair $(\alpha'_1,\alpha'_2)$ is representable in $\cB_{1,q_0X}$
by a morphism $\lambda':E\to Y'_1$ in $\cB_{1,q_0X}$. Set
$\lambda:=\beta_1\circ\lambda'$; it follows that
$\alpha_1\circ\lambda=\alpha_2\circ\lambda$. Moreover, for
$i=1,2$ we may find a morphism $\phi'_i:X\to vY'_i$ in
$\cB_{0,q_0X}$ such that $v(\beta_i)\circ\phi'_i=\phi_i$.
We compute :
$$
v(\beta_2\circ\alpha'_1)\circ\phi'_1=
v(\alpha_1\circ\beta_1)\circ\phi'_1=v(\alpha_1)\circ\phi_1=\phi_2=
v(\alpha_2)\circ\phi_1=v(\alpha_2\circ\beta_1)\circ\phi'_1=
v(\beta_2\circ\alpha'_2)\circ\phi'_1
$$
whence $v(\alpha'_1)\circ\phi'_1=v(\alpha'_2)\circ\phi'_1$, since
$v(\beta_2)$ is cartesian. However, $v(\lambda)$ is the equalizer
of the pair $(v(\alpha'_1),v(\alpha'_2))$, since $v_{q_0X}$ is left
exact; so finally $\phi'_1$ factors through a unique morphism
$\phi''_1:X\to vE$ and $v(\lambda'):vE\to vY'_1$. Thus, we get
a morphism $X/\lambda:(E,\phi''_1)\to(Y_1,\phi_1)$ in $X/v\cB_1$
such that $(X/\alpha_1)\circ(X/\lambda)=(X/\alpha_2)\circ(X/\lambda)$,
whence the contention.

(ii): In view of (i), it suffices to show that conditions (a)
and (b) of lemma \ref{lem_filtered-final}(i) hold for $\rho^X$.
To this aim, let $u_{q_0X}:\cA_{1,q_0X}\to\cA_{0,q_0X}$ be the
restriction of $u$ to the fibre categories, and $(Y,\phi)$
any object of $(\pi_0X/u\cA_1)^o$. Arguing as in the proof of
(i), we find $(Y',\phi')\in\Ob((\pi_0X/u_{q_0X}\cA_{1,q_0X})^o)$
and a morphism $(Y,\phi)\to(Y',\phi')$ in $(\pi_0X/u\cA_1)^o$.
But notice that $\rho^X$ restricts to an isomorphism of categories
$(X/v_{q_0X}\cB_{1,q_0X})^o\isom(\pi_0X/u_{q_0X}\cA_{1,q_0X})^o$.
Condition (a) is an immediate consequence. Lastly, let
$(Y,\phi)\in\Ob(X/v\cB_1)$, $(Z,\phi')\in\Ob(\pi_0X/u\cA_1)$,
and a pair of morphisms in $\pi_0X/u\cA_1$ :
$$
\pi_0X/\beta_i:\rho^X(Y,\phi)=(\pi_1Y,\pi_0(\phi))\to(Z,\phi')
\qquad
i=1,2.
$$
Hence, $\beta_1,\beta_2:\pi_1Y\to Z$ are morphisms of $\cB_1$
with $u(\beta_1)\circ\pi_0(\phi)=u(\beta_2)\circ\pi_0(\phi)$.
Arguing as in the proof of (i) we find a cartesian morphism
$\beta:Y'\to Y$ in $\cB_1$ and a morphism $\psi:X\to vY'$
in $\cB_{0,q_0X}$ with $v(\beta)\circ\psi=\phi$. Notice that
$p_1(\beta_1\circ\pi_1(\beta))=p_0(\phi')=
p_1(\beta_2\circ\pi_1(\beta))$. Then we may find
$Z''\in\Ob(\cA_{1,q_0X})$ and a cartesian
morphism $\gamma:Z''\to Z$ in $\cA_1$ such that
$p_1(\gamma)=p_0(\phi')$, and for $i=1,2$, a morphism
$\beta''_i:\pi_1Y'\to Z'$ in $\cA_{1,q_0X}$ such that
$\beta_1\circ\pi_1(\beta)=\gamma\circ\beta''_i$. Since
$\pi_1$ restricts to an isomorphism
$\cB_{1,q_0X}\isom\cA_{1,q_0X}$, there exist
$Z'\in\Ob(\cB_{1,q_0X})$ with $\pi_1Z'=Z''$, and unique
morphisms $\beta'_i:Y'\to Z'$ in $\cB_{1,q_0X}$ such that
$\pi_1(\beta'_i)=\beta''_i$ for $i=1,2$. We compute :
$$
\begin{aligned}
u(\gamma)\circ u\pi_1(\beta'_1)\circ\pi_0(\psi)=
u(\beta_1\circ\pi_1(\beta))\circ\pi_0(\psi)&\,=
u(\beta_1)\circ\pi_0(\phi) \\
&\,=u(\beta_2)\circ\pi_0(\phi) \\
&\,=u(\beta_2\circ\pi_1(\beta))\circ\pi_0(\psi) \\
&\,=u(\gamma)\circ u\pi_1(\beta'_2)\circ\pi_0(\psi).
\end{aligned}
$$
Since $u(\gamma)$ is cartesian and $\pi_0$ restricts to
an isomorphism $\cB_{0,q_0X}\isom\cA_{0,q_0X}$ , we conclude
that $v(\beta'_1)\circ\psi=v(\beta'_2)\circ\psi$. Let
$\lambda:E\to Y'$ be the equalizer of the pair
$(\beta'_1,\beta'_2)$ in $\cB_{1,q_0X}$; since $v_{q_0X}$ is
left exact, $v(\lambda):vE\to vY'$ represents the equalizer
of $(v(\beta'_1),v(\beta'_2))$ in $\cA_{1,q_0X}$, so finally
$\psi$ factor through $v(\lambda)$ and a unique morphism
$\mu:X\to vE$ in $\cA_{1,q_0X}$. We have
$(E,\mu)\in\Ob(X/v\cB_1)$, and a morphism
$X/(\beta\circ\lambda):(E,\mu)\to(Y,\phi)$ in $X/v\cB_1$
such that $\rho^X(X/(\beta\circ\lambda))$ equalizes
$\pi_0X/\beta_1$ and $\pi_0X/\beta_2$ in $(\pi_0X/u\cA_1)^o$.
This concludes the proof of condition (b).
\end{pfclaim}
\end{proof}

\subsection{Fibred topoi}\label{sec_fibred-topoi}
For every small category $I$, we consider next the $2$-category
$$
\sPsFun(I,\Topos)
$$
whose objects $T:I\to\Topos$ shall be called the {\em
fibred topoi over $I$}. The $1$-cells, that is the
pseudo-natural transformations $\omega:T\Rightarrow S$
shall be called {\em morphisms of fibred topoi over $I$},
and denoted by a simple arrow $T\to S$. The $2$-cells,
{\em i.e.} the modifications $\omega\leadsto\omega'$ in
$\sPsFun(I,\Topos)$ shall be called also
{\em transformations} of morphisms of fibred topoi over $I$.

\begin{remark}\label{rem_fibred-topoi}
(i)\ \
Recall that the definition of the $2$-category $\Topos$ involves
the choice of two universes $\sU,\sU'$ with $\sU\in\sU'$, so that
every object of $\Topos$ is both a $\sU$-topos and a small
$\sU'$-category (see definition \ref{def_morph-of-topoi}(iv)).
Then, consider any pseudo-functor
$$
F:I\to\sU'\tdu\bCat
$$
such that (a): $F_i$ is an $\sU$-topos for every $i\in\Ob(I)$,
and (b): for every morphism $\phi:i\to j$ of $I$ there exists
a morphism of topoi $f_\phi:F_i\to F_j$ with $F_\phi=f_{\phi*}$.
Then we get a fibred topos $T$ over $I$ such that
$T_i:=F_i$ for every $i\in\Ob(I)$ and $T_\phi:=f_\phi$ for
every morphism $\phi$ of $I$. This fibred topos depends on
the choice of the morphisms of topoi $f_\phi$, but any two
such choices yield isomorphic fibred topoi over $I$.

(ii)\ \
In particular, any {\em functor} $F:I\to\sU'\tdu\bCat$ fulfilling
(a) and (b) can be regarded as a fibred topos over $I$; in this
case, the coherence constraint $(\delta^T,\gamma^T)$ of the
resulting $T$ shall be given by identity automorphisms of functors :
$$
\delta^T_i=\one_{T_i}
\qquad
\gamma^T_{\phi,\psi}=\one_{T_\phi\circ T_\psi}
\qquad
\text{in $\sU'\tdu\bCat$}
$$
which are however {\em not necessarily} identity $2$-cells
of the $2$-category $\Topos$, since one does not
necessarily have $f^*_\psi\circ f^*_\phi=f^*_{\psi\circ\phi}$, nor
$f^*_{\one_i}=\one_{T_i}$. Hence, even when it is given by an
actual functor, $T$ will only in general be a (non-strict)
pseudo-functor, when regarded as a fibred topos.

(iii)\ \
From the pseudo-functors of remark \ref{rem_Topos-*}(ii)
we deduce two strict pseudo-functors :
$$
\sPsFun(I^o,\sU'\tdu\bCat)^o
\xleftarrow{[-]^*}\sPsFun(I,\Topos)
\xrightarrow{[-]_*}\sPsFun(I,\sU'\tdu\bCat).
$$

$\bullet$\ \
Namely, $[-]_*$ assigns to every fibred topos $T:I\to\Topos$ over
$I$ the composition
$$
T_*:=(-)_*\circ T:I\to\sU'\tdu\bCat
\qquad
i\mapsto T_i
\qquad
(i\xrightarrow{\phi}j)\mapsto(T_{\phi*}:T_i\to T_j)
$$
and to every morphism of fibred topoi $\omega:T\to S$ the
pseudo-natural transformation
$\omega_*:=(-)_**\omega:T_*\Rightarrow S_*$. To every
transformation $\Xi:\omega\leadsto\omega'$ it assigns
the modification
$$
\Xi_*:=(-)_*\circ\Xi:\omega_*\leadsto\omega'_*.
$$

$\bullet$\ \
Likewise, $[-]^*$ assigns to every fibred topos $T$ over $I$
the pseudo-functor
$$
T^*:=((-)^*\circ T)^o:I^o\to\sU'\tdu\bCat
\qquad
i^o\mapsto T^o_i
\qquad
(j^o\xrightarrow{\phi^o}i^o)\mapsto(T^{*o}_\phi:T_j^o\to T_i^o)
$$
and to every morphism of fibred topoi $\omega:T\to S$ the
pseudo-natural transformation
$\omega^*:=((-)^**\omega)^o:S^*\Rightarrow T^*$. To every
transformation $\Xi:\omega\leadsto\omega'$ it assigns the
modification
$$
\Xi^*:=((-)^*\circ\Xi)^o:\omega^*\leadsto\omega'^*.
$$
\end{remark}

\sset\subsubsection{}\label{subsec_two-fibrations}
There are several useful ways of attaching a fibration to any
fibred topos $T$ over $I$ : namely, with the notation of
\eqref{subsec_fib-from-pseudo} and remark
\ref{rem_fibred-topoi}(iii), we can consider the fibrations
$$
\pi:\cFib(T_*)\to I^o
\qquad\text{and}\qquad
\pi':\cFib(T^*)\to I.
$$
More precisely, we have strict pseudo-functors :
$$
\sU'\tdu\Fib(I)^o\xleftarrow{\cFib_I^o\circ[-]^*}\sPsFun(I,\Topos)
\xrightarrow{\cFib_{I^o}\circ[-]_*}\sU'\tdu\Fib(I^o).
$$

$\bullet$\ \
We have also a natural pseudo-functor from fibred topoi over
$I$ to fibred sites over $I$ : namely, the composition
$$
\underline\Can:\sPsFun(I,\Topos)\xrightarrow{\ \sPsFun(I,\Can)\ }
\sPsFun(I,\sU'\tdu\Site)\xrightarrow{\ \underline\cFib_I\ }
\sU'\tdu\fib.\Site(I)
$$
where $\Can:\Topos\to\Site$ is the pseudo-functor of
\eqref{subsec_T-and-Can}, and $\underline\cFib_I$ is the
$2$-equivalence of remark \ref{rem_fibred-site}(v).
Explicitly, $\underline\Can(T)$ is given as follows.
Since $I$ is a usual category, we have ${}^oI^o=I^o$,
so $T^*$ induces a pseudo-functor
$$
I^o\xrightarrow{{}^oT^*}{}^o(\sU'\tdu\bCat)
\xrightarrow{(-)^o}\sU'\tdu\bCat
\qquad
i^o\mapsto T_i
\qquad
(j^o\xrightarrow{\phi}i^o)\mapsto(T^*_\phi:T_j\to T_i)
$$
where $(-)^o$ denotes the strict isomorphism of example
\ref{ex_opposite-cat-as-pseudo-fun}. The fibration
$$
p:\underline\Can(T):=\cFib((-)^o\circ{}^oT^*)\to I
$$
carries a natural structure of fibred lex-site : namely,
the category $p^{-1}(i)$ is naturally identified with $T_i$ for
every $i\in\Ob(I)$, and inherits the latter's canonical topology;
also, $p$ admits a natural cleavage whose associated
pseudo-functor $\sc$ is identified with $(-)^o\circ{}^oT^*$,
so for every morphism $\phi:i\to j$ of $I$ the functor
$\sc_\phi:p^{-1}(j)\to p^{-1}(i)$ is in turn identified with the left
exact functor $T^*_\phi:T_j\to T_i$, which is continuous for the
canonical topologies (remark \ref{rem_choose-two-univs}(iii)).

$\bullet$\ \
Conversely, to every fibred site we may attach a fibred topos,
via the pseudo-functors
$$
\begin{aligned}
\underline\sT&\,:\sU'\tdu\fib.\Site(I)\xrightarrow{\ \underline\sc^\bullet\ }
\sPsFun(I,\sU'\tdu\Site)\xrightarrow{\ \sPsFun(I,\sT)\ }\sPsFun(I,\Topos) \\
\underline{\lex.\sT}&\,:\sU'\tdu\fib.\lex.\Site(I)
\xrightarrow{\ \lex.\underline\sc^\bullet\ }
\sPsFun(I,\sU'\tdu\lex.\Site)
\xrightarrow{\ \sPsFun(I,\lex.\sT)\ }\sPsFun(I,\Topos)
\end{aligned}
$$
where $\sT:\Site\to\Topos$ and $\lex.\sT:\lex.\Site\to\Topos$
are as in \eqref{subsec_T-and-Can}, and $\underline\sc^\bullet$
(resp. $\lex.\underline\sc^\bullet$) is a strict and strong
pseudo-inverse for $\underline\cFib_I$ (resp. for
$\lex.\underline\cFib_I$ : see remark \ref{rem_fibred-site}(vi)).

\begin{lemma}\label{lem_change-variance}
With the notation of \eqref{subsec_two-fibrations},
there exists a natural isomorphism of $I$-categories
$$
\cFib(T^*)\isom\cFib(T_*)^o.
$$
\end{lemma}
\begin{proof} For every morphism $\phi:i\to j$ of $I$,
let $(\eta_\phi,\eps_\phi)$ be the unit and counit for the
adjunction for the pair $(T^*_\phi,T_{\phi*})$ that defines
the morphism of topoi $T_\phi:T_i\to T_j$. This adjunction
assigns to every morphism $f:T^*_\phi Y\to X$ (for any
$X\in\Ob(T_i)$ and $Y\in\Ob(T_j)$) an adjoint morphism
$$
f^\dagger:=T_{\phi*}(f)\circ\eta_{\phi,Y}:Y\to T_{\phi*}X.
$$
Notice that $\cFib(T^*)$ and $\cFib(T_*)$ have the same
set of objects, and the morphisms of $\cFib(T^*)$ (resp.
$\cFib(T_*)$) are the pairs $(\phi,f)$ where $\phi:i\to j$
is a morphism of $I$ and $f:T^*_\phi Y\to X$ (resp.
$f:Y\to T_{\phi*}X$) is a morphism of $T_i$ (resp. of $T_j$).
We claim that the rules
$$
(i,X)\mapsto(i,X)
\qquad
(\phi,f)\mapsto(\phi,f^\dagger)
$$
yield the sought isomorphism of categories. Indeed, denote
also by $(\delta,\gamma)$ the coherence constraint of
the pseudo-functor $T^*$; hence
$$
\gamma_{(\phi,\psi)}:T^*_{\psi\phi}\isom T^*_\phi T^*_\psi
\qquad\text{and}\qquad
\delta_i:T^*_{\one_i}\isom\one_{T_i}
$$
are isomorphisms of functors, for every $i\in\Ob(I)$ and
every pair of morphisms $\phi:i\to j$ and $\psi:j\to k$
in $I$, and according to remark \ref{rem_adjoint-transf}(ii),
there follow adjoint isomorphisms
$$
\gamma^\dagger_{(\phi,\psi)}:T_{\psi*}T_{\phi*}\isom T_{\psi\circ\phi*}
\qquad\text{and}\qquad
\delta^\dagger_i:\one_{T_i}\isom T_{\one_i*}.
$$
With this notation, we need to check the identity
$$
(f\circ T^*_\phi g\circ\gamma_{(\phi,\psi),Z})^\dagger=
\gamma^\dagger_{(\phi,\psi),X}\circ T_{\psi*}(f^\dagger)\circ g^\dagger
$$
for every pair of morphisms $(\phi,f:T^*_\phi Y\to X)$ and
$(\psi,g:T^*_\psi Z\to Y)$ of $\cFib(T^*)$. However, on
the one hand we have
$$
(f\circ T^*_\phi(g)\circ\gamma_{(\phi,\psi),Z})^\dagger=
T_{\psi\phi*}(f\circ T^*_\phi g\circ\gamma_{(\phi,\psi),Z})\circ
\eta_{\psi\phi,Z}
$$
and on the other hand :
$$
\begin{aligned}
\gamma^\dagger_{(\phi,\psi),X}\circ T_{\psi*}(f^\dagger)\circ g^\dagger
=&\,\gamma^\dagger_{(\phi,\psi),X}\circ T_{\psi*}T_{\phi*}(f)\circ
T_{\psi*}(\eta_{\phi,Y})\circ T_{\psi*}(g)\circ\eta_{\psi,Z} \\
=&\,T_{\psi\phi*}(f)\circ\gamma^\dagger_{(\phi,\psi),T^*_\phi Y}\circ
T_{\psi*}T_{\phi*}T^*_\phi(g)\circ T_{\psi*}(\eta_{\phi,T^*_\psi Z})
\circ\eta_{\psi,Z} \\
=&\,T_{\psi\phi*}(f)\circ\gamma^\dagger_{(\phi,\psi),T^*_\phi Y}\circ
T_{\psi*}T_{\phi*}T^*_\phi(g)\circ\eta_{(\phi,\psi),Z}
\end{aligned}
$$
where $(\eta_{(\phi,\psi)},\eps_{(\phi,\psi)})$ denotes the
unit and counit of the adjunction that defines the
composition $T_\psi\circ T_\phi:T_i\to T_k$, by virtue
of remark \ref{rem_adjoint-transf}(i).
Thus, we are reduced to checking that :
$$
T_{\psi\phi*}(T^*_\phi g\circ\gamma_{(\phi,\psi),Z})\circ
\eta_{\psi\phi,Z}=
\gamma^\dagger_{(\phi,\psi),T^*_\phi Y}\circ T_{\psi*}T_{\phi*}T^*_\phi(g)
\circ\eta_{(\phi,\psi),Z}.
$$
But we have
$\gamma^\dagger_{(\phi,\psi),T^*_\phi Y}\circ T_{\psi*}T_{\phi*}T^*_\phi(g)=
T_{\psi\phi*}T^*_\phi(g)\circ\gamma^\dagger_{(\phi,\psi),T^*_\phi T^*_\psi Z}$,
so we are further reduced to showing the identity :
$$
T_{\psi\phi*}(\gamma_{(\phi,\psi),Z})\circ\eta_{\psi\phi,Z}=
\gamma^\dagger_{(\phi,\psi),T^*_\phi T^*_\psi Z}\circ\eta_{(\phi,\psi),Z}.
$$
But using the explicit expressions of remark
\ref{rem_adjoint-transf}(ii) we see that :
$$
\gamma^\dagger_{(\phi,\psi)}*T^*_\phi T^*_\psi=
(T_{\psi\phi*}*\eps_{(\phi,\psi)}*T^*_\phi T^*_\psi)\!\odot\!
(T_{\psi\phi*}*\gamma_{(\phi,\psi)}*T_{\psi*}T_{\phi*}T^*_\phi T^*_\psi)
\!\odot\!(\eta_{\psi\phi}*T_{\psi*}T_{\phi*}T^*_\phi T^*_\psi)
$$
and on the other hand
$$
(\eta_{\psi\phi}*T_{\psi*}T_{\phi*}T^*_\phi T^*_\psi)\odot\eta_{(\phi,\psi)}
=(T_{\psi\phi*}T^*_{\psi\phi}*\eta_{(\phi,\psi)})\odot\eta_{\psi\phi}
$$
so it suffices to check that :
$$
\gamma_{(\phi,\psi)}=(\eps_{(\phi,\psi)}*T^*_\phi T^*_\psi)\odot
(\gamma_{(\phi,\psi)}*T_{\psi*}T_{\phi*}T^*_\phi T^*_\psi)\odot
(T^*_{\psi\phi}*\eta_{(\phi,\psi)}).
$$
The latter follows from the identity :
$$
(\gamma_{(\phi,\psi)}*T_{\psi*}T_{\phi*}T^*_\phi T^*_\psi)\odot
(T^*_{\psi\phi}*\eta_{(\phi,\psi)})=(T^*_\phi T^*_\psi*\eta_{(\phi,\psi)})
\odot\gamma_{(\phi,\psi)}
$$
together with the triangular identities for
$(\eta_{(\phi,\psi)},\eps_{(\phi,\psi)})$.
\end{proof}

\begin{definition}
(i)\ \
For any two categories $\cA,\cB$ and any functor $F:\cA\to\cB$
we denote
$$
\Sigma(\cA/\cB)
$$
the {\em category of sections of $F$}, which is the
subcategory of $\bFun(\cB,\cA)$ whose objects are
the functors $G:\cB\to\cA$ such that $F\circ G=\one_\cB$
and whose morphisms are the natural transformations
$\beta:G\Rightarrow G'$ with $F*\beta=i_\cB$ (where
$i_\cB$ is the identity automorphism of $\one_\cB$).

(ii)\ \
Let $T$ be a fibred topos over $I$, and $\pi,\pi'$ its
associated fibrations as in \eqref{subsec_two-fibrations};
we let :
$$
\sTop(T)_*:=\Sigma(\cFib(T_*)/I^o)
\qquad\text{and}\qquad
\sTop(T):=\Sigma(\cFib(T^*)/I)^o.
$$
By virtue of lemma \ref{lem_change-variance} and remark
\ref{rem_opposite-Fun}(i) we have a natural isomorphism
of categories :
\set\begin{equation}\label{eq_two-versions-of-tot-topos}
\sTop(T)\isom\sTop(T)_*.
\end{equation}
\end{definition}

\sset\subsubsection{}\label{subsec_use-oldstuff}
Let $\omega:T\to S$ be any morphism of fibred topoi
over $I$. From the pseudo-natural transformations
$\omega^*:S^*\Rightarrow T^*$ and
$\omega_*:T_*\Rightarrow S_*$ we get the functors
$$
\begin{aligned}
\cFib(\omega^*)&:\cFib(S^*)\to\cFib(T^*) \\
\cFib(\omega_*)&:\cFib(T_*)\to\cFib(S_*)
\end{aligned}
$$
(notation of \eqref{subsec_sd}); whence functors
$$
\begin{aligned}
\sTop(\omega)^*&:\sTop(S)\to\sTop(T)
\qquad
E\mapsto\cFib(\omega^*)\circ E \\
\sTop(\omega)_*&:\sTop(T)_*\to\sTop(S)_*
\qquad
E\mapsto\cFib(\omega_*)\circ E.
\end{aligned}
$$
By virtue of \eqref{eq_two-versions-of-tot-topos} we may
then identify $\sTop(\omega)_*$ with a functor which we shall
also denote
$$
\sTop(\omega)_*:\sTop(T)\to\sTop(S).
$$
Moreover, $\omega$ also induces a morphism of fibred sites
$$
\underline\Can(\omega):=\cFib((-)^o*{}^o\omega^*):
\underline\Can(T)\to\underline\Can(S).
$$
Let $\omega:T\to S$ and $\nu:S\to U$ be two morphisms of
fibred topoi over $I$. A simple inspection of the definition
shows that
\set\begin{equation}\label{eq_compose-Tops}
\begin{aligned}
\sTop(\omega)^*\circ\sTop(\nu)^*&\,=\sTop(\nu\circ\omega)^* \\
\sTop(\nu)_*\circ\sTop(\omega)_*&\,=\sTop(\nu\circ\omega)_* \\
\underline\Can(\nu)\circ\underline\Can(\omega)&\,=
\underline\Can(\nu\circ\omega).
\end{aligned}
\end{equation}

\sset\subsubsection{}
Let $\omega,\omega':T\to S$ be two morphisms of fibred topoi
over $I$, and $\Xi:\omega\leadsto\omega'$ a transformation.
Then $\Xi$ induces natural transformations
$$
\cFib(\Xi_*):\cFib(\omega_*)\Rightarrow\cFib(\omega'_*)
\qquad
\cFib(\Xi^*):\cFib(\omega^*)\Rightarrow\cFib(\omega'^*)
$$
whence natural transformations :
$$
\begin{aligned}
\sTop(\Xi)^*&:\sTop(\omega')^*\Rightarrow\sTop(\omega)^*
\qquad
E\mapsto\cFib(\Xi^*)*E \\
\sTop(\Xi)_*&:\sTop(\omega)_*\Rightarrow\sTop(\omega')_*
\qquad
E\mapsto\cFib(\Xi_*)*E.
\end{aligned}
$$
Moreover, $\Xi$ induces a natural transformation of morphisms
of fibred sites :
$$
\underline\Can(\Xi):=\cFib((-)^o\circ{}^o\Xi^*):
\underline\Can(\omega')\to\underline\Can(\omega).
$$

\begin{remark}\label{rem_use-oldstuff}
(i)\ \
Let $T$ be a fibred topos over $I$ with coherence constraint
$(\delta^T,\gamma^T)$, and $T^u:I\to\Topos$ the unital
pseudo-functor associated with $T$ (proposition
\ref{prop_towards-2-yoneda}); by \eqref{eq_compose-Tops}
and remark \ref{rem_isom-in-psfun}(i), the pseudo-natural
isomorphism $\alpha^T:T\isom T^u$ induces an isomorphism of
categories $\sTop(\alpha^T)^*:\sTop(T^u)\isom\sTop(T)$. Then
the following description of $\sTop(T^u)$ --  obtained by
direct inspection -- also applies to $\sTop(T)$, up to
natural isomorphism of categories:
\begin{itemize}
\item
the objects of $\sTop(T^u)$ are the systems
$E_\bullet:=((E_i,E_\phi)~|~i\in\Ob(I),\phi\in\rMorph(I))$
with $E_i\in\Ob(T_i)$ for every $i\in\Ob(I)$ and where
$$
E_\phi:T^*_\phi E_j\to E_i
\qquad
\text{for every $\phi:i\to j$ in $I$}
$$
is a morphism of $T_i$, such that the following diagram
commutes :
$$
\xymatrix{ T_\psi^*T_\phi^*E_k \ar[rr]^-{T_\psi^*E_\phi} & &
T^*_\psi E_j \ar[d]^{E_\phi} \\
T_{\psi\circ\phi}^*E_k \ar[rr]^-{E_{\psi\circ\phi}}
\ar[u]^{(\gamma^{T^u}_{(\phi,\psi)})^\dagger_{E_k}} & & E_i
}$$
for every pair of morphisms
$i\xrightarrow{\phi}j\xrightarrow{\psi}k$ of $I$, and
$E_{\one_i}=\one_{E_i}$ for every $i\in\Ob(I)$.
\item
the morphisms $f_\bullet:E_\bullet\to F_\bullet$ are the systems
$(f_i~|~i\in\Ob(I))$ where $f_i:E_i\to F_i$ is a morphism
of $T_i$ for every $i\in\Ob(I)$, such that the following
diagram commutes :
$$
{\diagram T^*_\phi E_j \ar[rr]^-{T_\phi^*f_j} \ar[d]_{E_\phi} & &
T^*_\phi F_j \ar[d]^{F_\phi} \\
E_i \ar[rr]^-{f_i} & & F_i
\enddiagram}
\qquad
\text{for every morphism $\phi:i\to j$ of $I$}.
$$
\item
the composition of morphisms $f_\bullet:E_\bullet\to F_\bullet$
and $g_\bullet:F_\bullet\to G_\bullet$ is defined by the obvious
rule : $(g_\bullet\circ f_\bullet)_i:=g_i\circ f_i$ for every
$i\in\Ob(I)$.
\end{itemize}

(ii)\ \
Likewise, notice that for every unital fibred topos $T$
over $I$, the objects of $\sTop(T)_*$ are the systems
$F_\bullet:=(F_i,F_\phi~|~i\in\Ob(I),\ \phi\in\rMorph(I))$
with $F_i\in\Ob(T_i)$ for every $i\in\Ob(I)$ and
$F_\phi:F_j\to T_{\phi*}F_i$ a morphism of $T_j$, for every
morphism $\phi:i\to j$ of $I$. This datum is required to
satisfy the identity
$$
\gamma^T_{(\phi,\psi),F_i}\circ(T_{\phi*}F_\phi)\circ F_\psi=F_{\psi\circ\phi}
$$
for every pair $i\xrightarrow{\phi}j\xrightarrow{\psi}k$ of
morphisms of $I$, and $F_{\one_i}=\one_{F_i}$ for every
$i\in\Ob(I)$, where, as usual, $\gamma^T$ denotes
the coherence constraint of $T$. The morphisms
$f_\bullet:F_\bullet\to F'_\bullet$ are the systems
$(f_i:F_i\to F'_i~|~i\in\Ob(I))$ where $f_i$ is a morphism of
$T_i$ for every $i\in\Ob(I)$ and
$$
F'_\phi\circ f_j=(T_{\phi*}f_i)\circ F_\phi
\qquad
\text{for every morphism $\phi:i\to j$ of $I$}.
$$

(iii)\ \
To every morphism $\omega:T\to S$ of fibred topoi over $I$ there
corresponds a morphism of unital fibred topoi $\omega^u:T^u\to S^u$
(see remark \ref{rem_unital}(i)); thus, in order to describe
$\sTop(\omega)^*$, we may assume without loss of generality that
$T$ and $S$ are unital. In this case, we get the following explicit
description. By definition, the coherence constraint of $\omega$
assigns to every morphism $\phi:i\to j$ of $I$ a natural
transformation
$$
\xymatrix{
S_j \ar[r]^-{\omega^*_j} \ar[d]_{S^*_\phi}
\drtwocell\omit{_\ \ \ \ \tau_\phi^{\omega\dagger}} &
T_j \ar[d]^{T^*_\phi} \\
S_i \ar[r]_-{\omega^*_i} & T_i.
}$$
Then $\sTop(\omega)^*$ is given by the rule :
$E_\bullet\mapsto
((\omega^*_i(E_i),E^\omega_\phi)~|~i\in\Ob(I),\phi\in\rMorph(I))$,
with
$$
E^\omega_\phi:
T^*_\phi\circ\omega^*_j(E_j)
\xrightarrow{\ \tau^{\omega\dagger}_{\phi,E_j}\ }
\omega^*_i\circ S^*_\phi(E_j)
\xrightarrow{\ \omega^*_i(E_\phi)\ }\omega^*_i(E_i)
$$
and by assigning to every morphism $g_\bullet:E_\bullet\to E'_\bullet$
of $\sTop(S)$ the morphism
$$
(\omega^*_i(g_i):
\omega^*_i(E_i)\to\omega^*_i(E'_i)~|~i\in\Ob(I)).
$$
\end{remark}

\begin{proposition}\label{prop_top-adjoint-pair}
The functor $\sTop(\omega)^*$ is left adjoint to $\sTop(\omega)_*$.
\end{proposition}
\begin{proof} Arguing as in remark \ref{rem_use-oldstuff}(i),
we reduce easily to the case where $T$ and $S$ are unital. In
light of the identification \eqref{eq_two-versions-of-tot-topos},
it suffices to exhibit an adjunction for the pair of functors
$$
\xymatrix{ \sTop(S)_* \ar@<.5ex>[rr]^-{\Omega^*} & & \sTop(T)
\ar@<.5ex>[ll]^-{\Omega_*}
}$$
resulting from the pair $(\sTop(\omega)^*,\sTop(\omega)_*)$
via these identifications. Explicitly, $\Omega^*$ assigns
to every object $F_\bullet$ of $\sTop(S)_*$ (see remark
\ref{rem_use-oldstuff}(ii)) the following object of $\sTop(T)$ :
$$
(\omega^*_i(F_i),F^\Omega_\phi:=
\omega^*_i(F_\phi^\dagger)\circ(\tau^\omega_\phi)^\dagger_{F_j}~|~
i\in\Ob(I),\ \phi\in\rMorph(I))
$$
where $F^\dagger_\phi:T^*_\phi F_j\to F_i$ is the morphism of
$T_i$ corresponding to $F_\phi$ under the adjunction of $T_\phi$,
{\em i.e.} $F^\dagger_\phi=\eps^{T_\phi}_{E_i}\circ(T^*_\phi F_\phi)$,
where $\eps^{T_\phi}$ is the counit of the adjunction of $T_\phi$.
For every morphism $f_\bullet:F_\bullet\to F'_\bullet$ of $\sTop(S)_*$
we have $\Omega^*(f_\bullet)=(\omega^*_if_i~|~i\in\Ob(I))$.

Likewise, $\Omega_*$ assigns to every object $E_\bullet$ of
$\sTop(T)$ the following object of $\sTop(S)_*$ :
$$
(\omega_{i*}(E_i),E^\Omega_\phi:=
(\tau^\omega_\phi)^{-1}_{E_i}\circ\omega_{j*}(E^\dagger_\phi)~|~
i\in\Ob(I),\ \phi\in\rMorph(I))
$$
where $E^\dagger_\phi:E_j\to T_{\phi*}E_i$ is the morphism of $T_j$
corresponding to $E_\phi$ under the adjunction of $T_\phi$,
{\em i.e.} $E^\dagger_\phi=(T_{\phi*}E_\phi)\circ\eta^{T_\phi}_{E_j}$,
where $\eta^{T_\phi}$ is the unit of $T_\phi$. For every morphism
$g_\bullet:E_\bullet\to E'_\bullet$ of $\sTop(T)$ we have
$\Omega_*(g_\bullet)=(\omega_{i*}g_i~|~i\in\Ob(I))$.

Now, let $\beta:\Omega^*(F_\bullet)\to E_\bullet$ be a morphism
of $\sTop(T)$; thus, $\beta$ is a system of morphisms
$(\beta_i:\omega^*_i(F_i)\to E_i~|~i\in\Ob(I))$ such that
the following diagram commutes :
\set\begin{equation}\label{eq_ninja-turtle}
{\diagram T^*_\phi\omega^*_jF_j \ar[rr]^-{T^*_\phi\beta_j}
\ar[d]_{F^\Omega_\phi} & & T^*_\phi E_j \ar[d]^{E_\phi} \\
\omega^*_iF_i \ar[rr]^-{\beta_i} & & E_i
\enddiagram}
\qquad
\text{for every morphism $i\xrightarrow{\phi}j$ of $I$}.
\end{equation}
Our candidate adjunction assigns to $\beta$ the system of morphisms
$$
\beta^\dagger:=(\beta^\dagger_i:=
\omega_{i*}(\beta_i)\circ\eta^{\omega_i}_{F_i}:
F_i\to\omega_{i*}E_i~|~i\in\Ob(I))
$$
where $\eta^{\omega_i}$ is the unit of $\omega_i:T_i\to S_i$,
so $\beta^\dagger_i$ is the morphism corresponding to $\beta_i$
under the adjunction of $\omega_i$. Thus, we need to show the
commutativity of the diagram :
$$
{\diagram F_j \ar[rr]^-{\beta^\dagger_j} \ar[d]_{F_\phi} & &
\omega_{j*}E_j \ar[d]^{E^\Omega_\phi} \\
S_{\phi*}F_i \ar[rr]^-{S_{\phi*}(\beta^\dagger_i)} & &
S_{\phi*}\omega_{i*}E_i
\enddiagram}
\qquad
\text{for every morphism $i\xrightarrow{\phi}j$ of $I$}.
$$
Set $R:=\omega_j\circ T_\phi$, and let $\eta^R$ and
$\eps^R$ be the unit and counit of $R$; we consider the diagram :
$$
\xymatrix@C+40pt{ F_j \ar[r]^-{\eta^R_{F_j}} \ar[d]_{\beta^\dagger_j} &
R_*R^*F_j \ar[r]^-{R_*(F^\Omega_\phi)} \ar[d]_{R_*T^*_\phi(\beta_j)} &
R_*\omega^*_iF_i \ar[r]^-{{(\tau^\omega_\phi)_{\omega^*_iF_i}^{-1}}}
\ar[d]^{R_*(\beta_i)} & S_{\phi*}\omega_{i*}\omega^*_iF_i
\ar[d]^{S_{\phi*}\omega_{i*}(\beta_i)} \\
\omega_{j*}E_j \ar[r]_-{\omega_{j*}(\eta^{T_\phi}_{E_j})} &
R_*T^*_\phi E_j \ar[r]_-{R_*(E_\phi)} &
R_*E_i \ar[r]_-{(\tau^\omega_\phi)_{E_i}^{-1}} &
S_{\phi*}\omega_{i*}E_i.
}$$
Recalling that $\eta^R_{F_j}=
\omega_{j*}(\eta^{T_\phi}_{\omega^*_jF_j})\circ\eta^{\omega_j}_{F_j}$,
we see that the left square commutes, by naturality of
$\eta^{T_\phi}$. The commutativity of the central square follows
trivially from that of \eqref{eq_ninja-turtle}, and that
of the right square follows from the naturality of
$\tau^\omega_\phi$. Notice now that the composition of the
three bottom horizontal arrows equals $E^\Omega_\phi$.
We are therefore reduced to checking the identity :
$$
(\tau^\omega_\phi)_{\omega^*_iF_i}^{-1}\circ
R_*(F^\Omega_\phi)\circ\eta^R_{F_j}
=S_{\phi*}(\eta^{\omega_i}_{F_i})\circ F_\phi.
$$
We compute :
$$
\begin{aligned}
R_*(F^\Omega_\phi)\circ\eta^R_{F_j}&\,=
R_*(\omega^*_i(F^\dagger_\phi)\circ
\tau^{\omega\dagger}_{\phi,F_j})\circ\eta^R_{F_j} \\
&\,=R_*\omega^*_i(F^\dagger_\phi)\circ
R_*(\eps^R_{\omega^*_iS^*_\phi F_j})\circ
R_*R^*(\tau^\omega_{\phi,\omega^*_iS^*_\phi F_j}\circ
\eta^{S_\phi\circ\omega_i}_{F_j})\circ\eta^R_{F_j} \\
&\,=R_*\omega^*_i(F^\dagger_\phi)\circ
R_*(\eps^R_{\omega^*_iS^*_\phi F_j})\circ
\eta^R_{R_*\omega^*_iS^*_\phi F_j}\circ\tau^\omega_{\phi,\omega^*_iS^*_\phi F_j}
\circ\eta^{S_\phi\circ\omega_i}_{F_j} \\
&\,=R_*\omega^*_i(F^\dagger_\phi)\circ
\tau^\omega_{\phi,\omega^*_iS^*_\phi F_j}
\circ\eta^{S_\phi\circ\omega_i}_{F_j} \\
&\,=\tau^\omega_{\phi,\omega^*_iF_i}\circ
S_{\phi*}\omega_{i*}\omega^*_i(F^\dagger_\phi)\circ
\eta^{S_\phi\circ\omega_i}_{F_j} \\
&\,=\tau^\omega_{\phi,\omega^*_iF_i}\circ
S_{\phi*}\omega_{i*}\omega^*_i(\eps^{S_\phi}_{E_i})\circ
S_{\phi*}\omega_{i*}\omega^*_i S^*_\phi(F_\phi)\circ
\eta^{S_\phi\circ\omega_i}_{F_j} \\
&\,=\tau^\omega_{\phi,\omega^*_iF_i}\circ
S_{\phi*}\omega_{i*}\omega^*_i(\eps^{S_\phi}_{E_i})\circ
\eta^{S_\phi\circ\omega_i}_{S_{\phi*}F_i}\circ F_\phi.
\end{aligned}
$$
So we are further reduced to showing that :
$$
S_{\phi*}\omega_{i*}\omega^*_i(\eps^{S_\phi}_{E_i})\circ
\eta^{S_\phi\circ\omega_i}_{S_{\phi*}F_i}=S_{\phi*}(\eta^{\omega_i}_{F_i}).
$$
But we have
$$
\begin{aligned}
S_{\phi*}\omega_{i*}\omega^*_i(\eps^{S_\phi}_{E_i})\circ
\eta^{S_\phi\circ\omega_i}_{T_{\phi*}F_i}&\,=
S_{\phi*}\omega_{i*}\omega^*_i(\eps^{S_\phi}_{E_i})\circ
S_{\phi*}(\eta^{\omega_i}_{S_\phi^*S_{\phi*}F_i})\circ
\eta^{S_\phi}_{S_{\phi*}F_i} \\
&\,=S_{\phi*}(\omega_{i*}\omega^*_i(\eps^{S_\phi}_{E_i})\circ
\eta^{\omega_i}_{S_\phi^*S_{\phi*}F_i})\circ\eta^{S_\phi}_{S_{\phi*}F_i} \\
&\,=S_{\phi*}(\eta^{\omega_i}_{F_i})\circ S_{\phi*}(\eps^{S_\phi}_{F_i})
\circ\eta^{S_\phi}_{S_{\phi*}F_i}
\end{aligned}
$$
so it suffices to invoke the triangular identities for
the pair $(\eta^{S_\phi},\eps^{S_\phi})$ to conclude.

Likewise one shows that for every morphism
$\lambda:F_\bullet\to\Omega_*E_\bullet$ in $\sTop(S)_*$, the
system
$$
(\lambda^\dagger_i:=\eps^{\omega_i}_{E_i}\circ\omega^*_i(\lambda_i)
:\omega^*_iF_i\to E_i~|~i\in\Ob(I))
$$
is a morphism $\Omega^*F_\bullet\to E_\bullet$ in $\sTop(T)$ :
the verification shall be left to the reader. In view of
remark \ref{rem_adjoint-transf}(ii), it is then clear that
these two rules yield mutually inverse bijections between
$\Hom_{\sTop(T)}(\Omega^*(F_\bullet),E_\bullet)$ and
$\Hom_{\sTop(S)_*}(F_\bullet,\Omega_*(E_\bullet))$. Lastly, the
naturality of the rule
$(\beta:\Omega^*(F_\bullet)\to E_\bullet)\mapsto\beta^\dagger$
with respect to both $F_\bullet$ and $E_\bullet$ follows directly
from the same property for the adjunction of each morphism
$\omega_i$ (details left to the reader).
\end{proof}

\begin{remark}\label{rem_top-adjoint-pair}
In the situation of proposition \ref{prop_top-adjoint-pair},
denote by $\lambda_S:\sTop(S)\isom\sTop(S)_*$ the natural
isomorphism of \eqref{eq_two-versions-of-tot-topos}. Then
we have $\sTop(\omega)_*\circ\sTop(\omega)^*=
\lambda^{-1}_S\circ\Omega_*\circ\Omega^*\circ\lambda_S$,
and $\sTop(\omega)^*\circ\sTop(\omega)_*=\Omega^*\circ\Omega_*$.
It follows that from the proof of proposition
\ref{prop_top-adjoint-pair} we can extract an adjunction
for the pair $(\sTop(\omega)^*,\sTop(\omega)_*)$, whose
unit $\eta^{\sTop(\omega)}$ and counit $\eps^{\sTop(\omega)}$
are given explicitly by the rules :
$$
F_\bullet\mapsto(\eta^{\omega_i}_{F_i}~|~i\in\Ob(I))
\qquad
E_\bullet\mapsto(\eps^{\omega_i}_{E_i}~|~i\in\Ob(I))
$$
where $(\eta^{\omega_i},\eps^{\omega_i})$ are the unit and
counit for the adjoint pair $(\omega_i^*,\omega_{i*})$,
for every $i\in\Ob(I)$.
\end{remark}

\sset\subsubsection{}\label{subsec_change-of-indexing-cat}
Let $\pi:\cA\to I$ and $\rho:J\to I$ be any two functors;
we have an obvious functor
$$
\Sigma(\cA/\rho)^*:\Sigma(\cA/I)\to\Sigma(J\times_I\cA/J)
\qquad
(G:I\to\cA)\mapsto J\times_IG
$$
that assigns to every morphism $\alpha:G\Rightarrow G'$ in
$\Sigma(\cA/I)$ the natural transformation
$J\times_I\alpha:J\times_IG\Rightarrow J\times_IG'$.
Especially, if $I$ and $J$ are small categories, and $T$ is
any fibred topos over $I$, recalling the natural identification
$\cFib((T\circ\rho)^*)\isom J\times_I\cFib(T^*)$ of remark
\ref{rem_added-little-extra}(ii), we deduce a functor
$$
\sTop(T/\rho)^*:=\Sigma(\cFib(T^*)/\rho)^{*o}:
\sTop(T)\to\sTop(T\circ\rho).
$$
Explicitly, $\sTop(T/\rho)^*$ assigns to every $E_\bullet$
as in remark \ref{rem_use-oldstuff}(i) the datum
$(E_{\rho(j)},E_{\rho(\phi)}~|~j\in\Ob(J),\ \phi\in\rMorph(J))$,
and to every morphism $f_\bullet:E_\bullet\to F_\bullet$ of
$\sTop(T)$, the morphism of $\sTop(T\circ\rho)$ given by the
system $(f_{\rho(j)}:E_{\rho(j)}\to F_{\rho(j)}~|~j\in\Ob(J))$.
We wish next to exhibit left and right adjoints for
$\sTop(T/\rho)^*$. To this aim, we observe more generally :

\begin{proposition}\label{adjoint-of-Sigma}
Let $\rho:J\to I$ be a functor between small categories, $\sU'$
a universe containing $\sU$, and $\sc:I^o\to\sU'\tdu\bCat$ any
pseudo-functor such that
\begin{enumerate}
\alphaenu
\item
The category $\sc_i$ is complete for every $i\in\Ob(I)$.
\item
The functor $\sc_f:\sc_i\to\sc_{i'}$ commutes with small
limits, for every morphism $i'\xrightarrow{f}i$ in $I$.
\end{enumerate}
Then $\Sigma(\cFib(\sc)/\rho)^*:
\Sigma(\cFib(\sc)/I)\to\Sigma(\cFib(\sc\circ\rho^o)/J)$
admits a right adjoint.
\end{proposition}
\begin{proof} By proposition \ref{prop_towards-2-yoneda}
we may assume that $\sc$ is unital. Now, set $\cA:=\cFib(\sc)$,
denote by $\pi:\cA\to I$ the projection, and let
$F_\bullet:J\to J\times_I\cA$ be any object of
$\Sigma(J\times_I\cA/J)$; hence $F_\bullet$ assigns
to every $j\in\Ob(J)$ an object $F_j\in\sc_{\rho(j)}$, and
to every morphism $\psi:j\to j'$ in $J$ a morphism
$F_\psi:F_j\to\sc_{\rho(\psi)}F_{j'}$ of $\sc_{\rho(j)}$. We set
$$
\Lambda_i(j,\phi):=\sc_\phi F_j
\qquad
\text{for every $i\in\Ob(I)$ and
$(j,\phi:i\to\rho(j))\in\Ob(i/\rho J)$}.
$$
For every morphism $(j,\phi)\xrightarrow{i/\psi}(j',\phi')$
of $i/\rho J$ we define $\Lambda_i(i/\psi):
\Lambda_i(j,\phi)\to\Lambda_i(j',\phi')$ as the composition
$$
\sc_\phi F_j\xrightarrow{\sc_\phi F_\psi}\sc_\phi\sc_{\rho(\psi)}F_{j'}
\xrightarrow{(\gamma^\sc_{\phi,\rho(\psi)})_{F_{j'}}}\sc_{\phi'} F_{j'}
$$
where $\gamma^\sc$ is the coherence constraint of $\sc$. Also, for
every morphism $f:i'\to i$ in $I$, let us set
$$
\Lambda_f(j,\phi):=(\gamma^\sc_{f,\phi})^{-1}_{F_j}:
\Lambda_{i'}(j,\phi\circ f)\to\sc_f(\Lambda_i(j,\phi))
\qquad
\text{for every $(j,\phi)\in\Ob(i/\rho J)$}.
$$

\begin{claim}\label{cl_tea-break}
(i)\ \
The rules : $(j,\phi)\mapsto\Lambda_i(j,\phi)$ and
$i/\psi\mapsto\Lambda(i/\psi)$ for every object $(j,\phi)$
and morphism $i/\psi$ of $i/\rho J$ define a functor
$$
\Lambda^{F_\bullet}_i:i/\rho J\to\pi^{-1}(i)
\qquad
\text{for every $i\in\Ob(I)$}.
$$

(ii)\ \
Let the functor $f/\rho J:i/\rho J\to i'/\rho J$ be as in
\eqref{subsec_fibreovercat}. The rule :
$(j,\phi)\mapsto\Lambda_f(j,\phi)$ defines a natural
transformation
$$
\Lambda^{F_\bullet}_f:\Lambda^{F_\bullet}_{i'}\circ(f/\rho J)
\Rightarrow\sc_f\circ\Lambda^{F_\bullet}_i.
$$
\end{claim}
\begin{pfclaim}(i): Since $\sc$ is unital, a simple inspection
shows that $\Lambda_i(\one_{(j,\phi)})=\one_{\Lambda(j,\phi)}$ for
every $(j,\phi)\in\Ob(\rho J/i)$. Next consider a pair
of morphisms $(j,\phi)\xrightarrow{i/\psi}(j',\phi')
\xrightarrow{i/\psi'}(j'',\phi'')$ of $i/\rho J$. We have
$$
\begin{aligned}
\Lambda_i(i/\psi')\circ\Lambda_i(i/\psi)&\,=
(\gamma^\sc_{\phi',\rho(\psi')})_{F_{j''}}\circ\sc_{\phi'}F_{\psi'}\circ
(\gamma^\sc_{\phi,\rho(\psi)})_{F_{j'}}\circ\sc_\phi F_\psi \\
&\,=(\gamma^\sc_{\phi',\rho(\psi')})_{F_{j''}}\circ
(\gamma^\sc_{\phi,\rho(\psi)}*\sc_{\rho(\psi')})_{F_{j'}}\circ
\sc_\phi\sc_{\rho(\psi)}(F_{\psi'})\circ\sc_\phi F_\psi \\
&\,=(\gamma^\sc_{\phi',\rho(\psi')})_{F_{j''}}\circ
(\gamma^\sc_{\phi,\rho(\psi)}*\sc_{\rho(\psi')})_{F_{j'}}\circ
\sc_\phi((\gamma^{\sc\ -1}_{\rho(\psi),\rho(\psi')})_{F_{j''}}\circ
F_{\psi'\circ\psi}) \\
&\,=(\gamma^\sc_{\phi,\rho(\psi'\circ\psi)})_{F_{j''}}\circ
\sc_\phi F_{\psi'\circ\psi} \\
&\,=\Lambda_i(i/\psi'\circ\psi)
\end{aligned}
$$
whence the contention.

(ii): It suffices to remark the commutativity of the diagram :
$$
\xymatrix{ \sc_{\phi\circ f}F_j \ar[rrr]^-{\sc_{\phi\circ f}F_\psi}
\ar[d]_{(\gamma^\sc_{f,\phi})^{-1}_{F_j}} & & &
\sc_{\phi\circ f}\sc_{\rho(\psi)}F_{j'}
\ar[d]^{(\gamma^\sc_{f,\phi})^{-1}_{\sc_{\rho(\psi)}F_{j'}}}
\ar[rrr]^-{(\gamma^\sc_{\phi\circ f,\rho(\psi)})_{F_{j'}}} & & &
\sc_{\rho(\psi)\circ\phi\circ f}F_{j'} \ar[d]^{(\gamma^\sc_{f,\phi'})^{-1}_{F_{j'}}} \\
\sc_f\sc_\phi F_j \ar[rrr]_-{\sc_f\sc_\phi F_\psi} & & &
\sc_f\sc_\phi\sc_{\rho(\psi)}F_{j'}
\ar[rrr]_-{\sc_f(\gamma^\sc_{\phi,\rho(\psi)})_{F_{j'}}} & & &
\sc_f\sc_{\rho(\psi)\circ\phi}F_{j'}
}$$
for every morphism $(i/\psi):(j,\phi)\to(j',\phi')$ of
$i/\rho J$.
\end{pfclaim}

In view of claim \ref{cl_tea-break}(i), for every $i\in\Ob(I)$
we choose $\lambda(F_\bullet,i)\in\Ob(\sc_i)$ representing the limit
of $\Lambda^{F_\bullet}_i$, and a universal cone
$\tau^{F_\bullet,i}:c_{\lambda(F_\bullet,i)}\Rightarrow\Lambda^{F_\bullet}_i$.
In view of claim \ref{cl_tea-break}(ii) there follows a cone
$$
\Upsilon^{F_\bullet}_f:=
\Lambda^{F_\bullet}_f\odot(\tau^{F_\bullet,i'}*(f/\rho J)):
c_{\lambda(F_\bullet,i')}\Rightarrow\sc_f\circ\Lambda^{F_\bullet}_i.
$$
On the other hand, notice that the cone
$\sc_f*\tau^{F_\bullet,i}:c_{\sc_f\lambda(F_\bullet,i)}\Rightarrow
\sc_f\circ\Lambda^{F_\bullet}_i$ is still universal, since
$\sc_f$ commutes with small limits; hence there exists a unique
morphism in $\sc_{i'}$
$$
\lambda(F_\bullet,f):\lambda(F_\bullet,i')\to\sc_f\lambda(F_\bullet,i)
\qquad\text{such that}\qquad
(\sc_f*\tau^{F_\bullet,i})\odot c_{\lambda(F_\bullet,f)}=
\Upsilon^{F_\bullet}_f.
$$

\begin{claim}\label{cl_amelie-jouit}
The rules : $i\mapsto\lambda(F_\bullet,i)$ and
$f\mapsto\lambda(F_\bullet,f)$ for every $i\in\Ob(I)$ and every
morphism $f$ of $I$ define a functor
$$
\lambda(F_\bullet):I\to\cA.
$$
\end{claim}
\begin{pfclaim} Since $\sc$ is unital, a simple inspection shows
that $\Lambda^{F_\bullet}_{\one_i}=\one_{\Lambda^{F_\bullet}_i}$ for every
$i\in\Ob(I)$, whence
$\lambda(F_\bullet,\one_i)=\one_{\lambda(F_\bullet,i)}$. Next, let
$i''\xrightarrow{f'}i'\xrightarrow{f}i$ be two morphisms of $I$;
we need to check that $\lambda(F_\bullet,f)\circ\lambda(F_\bullet,f')
=\lambda(F_\bullet,f\circ f')$, and by the universality of
$\sc_{f\circ f'}*\tau^{F_\bullet,i}$ it suffices to show that
$X:=(\sc_{f\circ f'}*\tau^{F_\bullet,i})\odot c_{\lambda(F_\bullet,f)\circ
\lambda(F_\bullet,f')}=(\sc_{f\circ f'}*\tau^{F_\bullet,i})\odot
c_{\lambda(F_\bullet,f\circ f')}$. We compute :
$$
\begin{aligned}
X&\,=(\sc_{f\circ f'}*\tau^{F_\bullet,i})\odot
c_{(\gamma^T_{f',f})_{\lambda(F_\bullet,i)}}\odot c_{\sc_{f'}\lambda(F_\bullet,f)}
\odot c_{\lambda(F_\bullet,f')} \\
&\,=(\gamma^\sc_{f,f'}*\Lambda^{F_\bullet}_i)\odot
(\sc_{f'}\sc_f*\tau^{F_\bullet,i})\odot c_{\sc_{f'}\lambda(F_\bullet,f)}
\odot c_{\lambda(F_\bullet,f')} \\
&\,=(\gamma^\sc_{f,f'}*\Lambda^{F_\bullet}_i)\odot
\sc_{f'}*(\Lambda^{F_\bullet}_f\odot(\tau^{F_\bullet,i'}*(f/\rho J)))
\odot c_{\lambda(F_\bullet,f')} \\
&\,=(\gamma^\sc_{f,f'}*\Lambda^{F_\bullet}_i)\odot
(\sc_{f'}*\Lambda^{F_\bullet}_f)\odot
((\Lambda^{F_\bullet}_{f'}\odot
(\tau^{F_\bullet,i''}*(f'/\rho J)))*(f/\rho J)) \\
&\,=(\gamma^\sc_{f,f'}*\Lambda^{F_\bullet}_i)\odot
(\sc_{f'}*\Lambda^{F_\bullet}_f)\odot
(\Lambda^{F_\bullet}_{f'}*(f/\rho J))\odot
(\tau^{F_\bullet,i''}*(f\circ f'/\rho J)).
\end{aligned}
$$
So we are reduced to checking that
$$
(\gamma^\sc_{f,f'}*\Lambda^{F_\bullet}_i)\odot
(\sc_{f'}*\Lambda^{F_\bullet}_f)\odot
(\Lambda^{F_\bullet}_{f'}*(f/\rho J))=\Lambda^{F_\bullet}_{f\circ f'}.
$$
But the latter follows directly from the coherence axioms
for $\gamma^\sc$.
\end{pfclaim}

It is clear that the functor $\lambda(F_\bullet)$ of claim
\ref{cl_amelie-jouit} is a section of the projection
$\cA\to I$. Next, let $\beta_\bullet:F_\bullet\to F'_\bullet$
be any morphism of $\Sigma(J\times_I\cA/J)$; we deduce
easily a natural transformation
$$
\Lambda^{\beta_\bullet}_i:\Lambda^{F_\bullet}_i\Rightarrow\Lambda^{F'_\bullet}_i
\qquad
(j,\phi)\mapsto(\sc_\phi\beta_j:\sc_\phi F_j\to\sc_\phi F'_j)
\qquad
\text{for every $i\in\Ob(I)$}
$$
such that
\set\begin{equation}\label{eq_natural}
(\sc_f\Lambda^{\beta_\bullet}_i)\odot\Lambda^{F_\bullet}_f=
\Lambda^{F'_\bullet}_f\odot(\Lambda^{\beta_\bullet}_{i'}*(f/\rho J))
\qquad
\text{for every morphism $f:i'\to i$ of $I$}
\end{equation}
whence a unique morphism in $\sc_i$
$$
\lambda(\beta_\bullet,i):\lambda(F_\bullet,i)\to\lambda(F'_{\!\bullet},i)
\qquad\text{such that}\qquad
\tau^{F'_\bullet,i}\odot c_{\lambda(\beta_\bullet,i)}=
\Lambda^{\beta_\bullet}_i\odot\tau^{F_\bullet,i}.
$$

\begin{claim} The rule $i\mapsto\lambda(\beta_\bullet,i)$
defines a natural transformation
$$
\lambda(\beta_\bullet):\lambda(F_\bullet)\Rightarrow\lambda(F'_\bullet).
$$
\end{claim}
\begin{pfclaim} Let $f:i'\to i$ be any morphism of $I$; as
usual, we reduce to checking the identity :
$X:=(\sc_f*\tau^{F'_\bullet,i})\odot c_{\lambda(F'_\bullet,f)}\odot
c_{\lambda(\beta_\bullet,i')}=Y:=(\sc_f*\tau^{F'_\bullet,i})\odot
(\sc_f*c_{\lambda(\beta_\bullet,i)})\odot c_{\lambda(F_\bullet,f)}$.
We compute :
$$
\begin{aligned}
X&\,=\Upsilon^{F'_\bullet}_f\odot c_{\lambda(\beta_\bullet,i')}=
\Lambda^{F'_\bullet}_f\odot(\Lambda^{\beta_\bullet}_{i'}\odot
\tau^{F_\bullet,i'})*(f/\rho J) \\
Y&\,=\sc_f(\Lambda^{\beta_\bullet}_i\odot\tau^{F_\bullet,i})\odot
c_{\lambda(F_\bullet,f)}=\sc_f\Lambda^{\beta_\bullet}_i\odot\Upsilon^{F_\bullet}_f
\end{aligned}
$$
so the assertion follows from \eqref{eq_natural}.
\end{pfclaim}

Let $\beta_\bullet:F_\bullet\to F'_\bullet$ and
$\beta'_\bullet:F'_\bullet\to F''_\bullet$ be two morphisms of
$\Sigma(J\times_I\cA/J)$; by a simple inspection
we see that $\Lambda^{\beta'_\bullet\circ\beta_\bullet}_i=
\Lambda^{\beta'_\bullet}_i\odot\Lambda^{\beta_\bullet}_i$ for every
$i\in\Ob(I)$, whence $\lambda(\beta'_\bullet\circ\beta_\bullet)=
\lambda(\beta'_\bullet)\odot\lambda(\beta_\bullet)$. It is also
easily seen that $\lambda(\one_{F_\bullet})=\one_{\lambda(F_\bullet)}$
for every object $F_\bullet$ of $\Sigma(J\times_I\cA/J)$,
so we have finally obtained a functor
$$
\Sigma(\cA/\rho)_*:\Sigma(J\times_I\cA/J)\to\Sigma(\cA/I)
\qquad
F_\bullet\mapsto\lambda(F_\bullet)
\qquad
(\beta_\bullet:F_\bullet\to F'_\bullet)\mapsto\lambda(\beta_\bullet).
$$
To see that $\Sigma(\cA/\rho)_*$ is the sought right adjoint,
consider any objects $E_\bullet$ of $\Sigma(\cA/I)$ and
$F_\bullet$ of $\Sigma(J\times_I\cA/J)$, and a morphism
$\beta_\bullet:E_\bullet\to\lambda(F_\bullet)$ in $\Sigma(\cA/I)$.
It is easily seen that the rule :
$$
j\mapsto\beta^\dagger_j:E_{\rho(j)}\xrightarrow{\beta_{\rho(j)}}
\lambda(F_\bullet,\rho(j))
\xrightarrow{\tau^{F_\bullet,\rho(j)}_{(j,\one_{\rho(j)})}}F_j
\qquad
\text{for every $j\in\Ob(J)$}
$$
defines a morphism $\beta^\dagger_\bullet:J\times_IE_\bullet\to F_\bullet$
in $\Sigma(J\times_I\cA/J)$ (details left to the reader).
Conversely, let $\alpha_\bullet:J\times_IE_\bullet\to F_\bullet$
be a morphism in $\Sigma(J\times_I\cA/J)$; for every
$i\in\Ob(I)$ and $(j,\phi)\in\Ob(i/\rho J)$ we let
$\tau^{\alpha_\bullet,i}_{(j,\phi)}:=\sc_\phi\alpha_j\circ E_\phi:
E_i\to\sc_\phi F_j$. It is easily seen that the rule :
$(j,\phi)\mapsto\tau^{\alpha_\bullet,i}_{(j,\phi)}$ defines a cone
$\tau^{\alpha_\bullet,i}:c_{E_i}\Rightarrow\Lambda^{F_\bullet}_i$
(details left to the reader); there follows a unique morphism
$$
\alpha^\dagger_i:E_i\to\lambda(F_\bullet,i)
\qquad\text{in $\sc_i$ such that}\qquad
\tau^{F_\bullet,i}\odot c_{\alpha^\dagger_i}=\tau^{\alpha_\bullet,i}.
$$

\begin{claim} The rule $i\mapsto\alpha^\dagger_i$ defines
a morphism $\alpha^\dagger_\bullet:E_\bullet\to\lambda(F_\bullet)$
in $\Sigma(\cA/I)$.
\end{claim}
\begin{pfclaim} As usual we reduce to checking that
$$
X:=(\sc_f*\tau^{F_\bullet,i})\odot c_{\sc_f\alpha^\dagger_i}\odot c_{E_f}=
Y:=(\sc_f*\tau^{F_\bullet,i})\odot c_{\lambda(F_\bullet,f)}\odot
c_{\alpha^\dagger_{i'}}
$$
for every morphism $f:i'\to i$ in $I$. However, we have :
$$
X=(\sc_f*\tau^{\alpha,i})\odot c_{E_f}
\qquad
Y=\Upsilon^{F_\bullet}_f\odot c_{\alpha^\dagger_{i'}}=
\Lambda^{F_\bullet}_f\odot(\tau^{\alpha,i'}*(f/\rho J)).
$$
So we come down to checking that
$$
(\sc_f\sc_\phi\alpha_j)\circ(\sc_fE_\phi)\circ E_f=
(\gamma^\sc_{f,\phi})^{-1}_{F_j}\circ\sc_{\phi\circ f}\alpha_j\circ
E_{\phi\circ f}
\qquad
\text{for every $(j,\phi)\in\Ob(i/\rho J)$}.
$$
To this aim, it suffices to notice that we have :
$(\gamma^\sc_{f,\phi})_{F_j}\circ(\sc_f\sc_\phi\alpha_j)=
\sc_{\phi\circ f}\alpha_j\circ(\gamma^\sc_{f,\phi})_{E_{\rho(j)}}$ and
$(\gamma^\sc_{f,\phi})_{E_{\rho(j)}}\circ(\sc_fE_\phi)\circ E_f=
E_{\phi\circ f}$.
\end{pfclaim}

\begin{claim}\label{cl_double-dagger}
For every morphism $\alpha_\bullet:J\times_IE_\bullet\to F_\bullet$
in $\Sigma(J\times_I\cA/J)$ and
$\beta_\bullet:E_\bullet\to\lambda(F_\bullet)$ in $\Sigma(\cA/I)$
we have $(\alpha^\dagger_\bullet)^\dagger=\alpha_\bullet$ and
$(\beta^\dagger_\bullet)^\dagger=\beta_\bullet$.
\end{claim}
\begin{pfclaim} The assertion concerning $\alpha_\bullet$
follows by direct inspection. Next, let us consider for
every $i\in\Ob(I)$ and every $(j,\phi)\in\Ob(i/\rho J)$
the diagram of morphisms in $\sc_i$
$$
\cD_{(j,\phi)}
\qquad : \qquad
{\diagram E_i \ar[rr]^-{E_\phi} \ar[rrd]_{\beta_i} & &
\sc_\phi E_{\rho(j)} \ar[rr]^-{\sc_\phi\beta_{\rho(j)}} & &
\sc_\phi\lambda(F_\bullet,\rho(j))
\ar[d]^{\sc_\phi\tau^{F_\bullet,\rho(j)}_{(j,\one_{\rho(j)})}} \\
& & \lambda(F_\bullet,i) \ar[rr]_-{\tau^{F_\bullet,i}_{(j,\phi)}}
\ar[rru]^-{\lambda(F_\bullet,\phi)} & & \sc_\phi F_j.
\enddiagram}$$
The assertion is equivalent to the commutativity of $\cD_{(j,\phi)}$
for every such $i$ and $(j,\phi)$. However, notice that
$(\Lambda^{F_\bullet}_\phi)_{(j,\one_{\rho(j)})}=\one_{\sc_\phi F_j}$,
since $\sc$ is unital; this implies that the lower triangular
subdiagram commutes. The commutativity of the upper triangular
subdiagram is clear.
\end{pfclaim}

Lastly, it is easily seen that the rule :
$\beta_\bullet\mapsto\beta^\dagger_\bullet$ is natural in both
$E_\bullet$ and $F_\bullet$ (details left to the reader); in
view of claim \ref{cl_double-dagger}, this rule then establishes
the sought adjunction between $\Sigma(\cA/\rho)^*$ and
$\Sigma(\cA/\rho)_*$.
\end{proof}

\begin{corollary}\label{cor_left-right-adj-Top-rho}
In the situation of \eqref{subsec_change-of-indexing-cat},
the functor $\sTop(T/\rho)^*$ admits both a left and a right
adjoint, denoted respectively :
$$
\sTop(T/\rho)_!:\sTop(T\circ\rho)\to\sTop(T)
\qquad\text{and}\qquad
\sTop(T/\rho)_*:\sTop(T\circ\rho)\to\sTop(T).
$$
\end{corollary}
\begin{proof} It is easily seen that the isomorphism of
categories \eqref{eq_two-versions-of-tot-topos} (and
the corresponding one for $T\circ\rho$) identifies the functor
$\sTop(T/\rho)^*$ with $\Sigma(\cFib(T_*)/\rho^o)^*$, hence
the existence of the right adjoint $\sTop(T/\rho)_*$ follows
from proposition \ref{adjoint-of-Sigma}. On the other hand,
by applying the same proposition to
$\Sigma(\cFib(T^*)/\rho)^*:\Sigma(\cFib(T^*)/I)\to
\Sigma(\cFib(T\circ\rho)^*/J)$ we see that the latter
admits as well a right adjoint, hence its opposite functor
$\sTop(T/\rho)^*$ admits a left adjoint.
\end{proof}

\begin{remark}\label{rem_fibrewise-lims-in-tot-Top}
(i)\ \
Let $\rho:J\to I$ be any functor between small categories
and $T$ any fibred topos over $I$. Corollary
\ref{cor_left-right-adj-Top-rho} implies that the functor
$\sTop(T/\rho)^*$ commutes with all small limits and all
small colimits (proposition \ref{prop_was-get-maddd}(iii,iv)).

(ii)\ \
Let $\bone$ be a final object of the category $\bCat$
({\em i.e.} a category with one object and one
morphism); for every $t\in\Ob(I)$ we have a unique functor
$\rho^t:\bone\to I$ that sends the unique object of $\bone$
to $t$. Clearly $\sTop(T\circ\rho^t)=T_t$ and the induced
functor
$$
\sTop(T/\rho^t)^*:\sTop(T)\to T_t
$$
assigns to every datum $E_\bullet$ as in remark
\ref{rem_use-oldstuff}(i) the object $E_t$ and to every
morphism $f_\bullet:E_\bullet\to F_\bullet$ of $\sTop(T)$ the
morphism $f_t:E_t\to F_t$ of $T_t$. The assertion that
$\sTop(T/\rho^t)^*$ commutes with limits and colimits for
every $t\in\Ob(I)$ then means that {\em the small limits
and colimits in the category $\sTop(T)$ are computed fibrewise}.

(iii)\ \
Moreover, the observation of (ii) determines the limits
and colimits in $\sTop(T)$ up to unique isomorphism. Indeed,
let $F:\Lambda\to\sTop(T)$ be any functor from a small
category $\Lambda$, and let $(L_\bullet,\tau)$ be a pair
consisting of an object $L_\bullet$ of $\sTop(T)$ representing
the colimit of $F$, and a universal cocone
$\tau:F\Rightarrow c_{L_\bullet}$. By (ii) we know that
$L_t\in\Ob(T_t)$ represents the colimit of
$F^t:=\sTop(T/\rho^t)^*\circ F:\Lambda\to T_t$ and
$\tau^t:=\sTop(T/\rho^t)^**\tau:F^t\Rightarrow c_{L_t}$ is
a universal cocone for every $t\in\Ob(I)$. Now, let
$\phi:s\to t$ be any morphism in $I$; we get a natural
transformation
$$
F^\phi:T^*_\phi\circ F^t\Rightarrow F^s
\qquad
\lambda\mapsto((F\lambda)_\phi:T^*_\phi(F\lambda)_t\to(F\lambda)_s).
$$
Since $T^*_\phi$ admits a right adjoint, the cocone
$T^*_\phi*\tau^t:T^*_\phi\circ F^t\Rightarrow c_{T^*_\phi L_t}$
is still universal (proposition \ref{prop_was-get-maddd}(iv)),
hence there exists a unique morphism
$$
L'_\phi:T^*_\phi L_t\to L_s
\qquad\text{in $T_s$ such that}\qquad
\tau^s\odot F^\phi=c_{L'_\phi}\odot(T^*_\phi*\tau^t).
$$
But then we must have $L'_\phi=L_\phi$ necessarily, so $L_\bullet$
is completely determined by a choice of universal cocones
$(\tau^t~|~t\in\Ob(I))$. Likewise we see that the limit of $F$
is completely determined by the choice of a system of pairs
$(M_t,\mu^t:F^t\Rightarrow c_{M_t}~|~t\in\Ob(I)$ with $M_t\in\Ob(T_t)$
representing the limit of $F^t$, and $\mu_t$ a universal cone,
for every $t\in\Ob(I)$.

(iv)\ \
By direct inspection, we see that the functor
$$
\sTop(T/\rho^t)_!:T_t\to\sTop(T)
$$
assigns to every $E\in\Ob(T_t)$ the datum
$(E_i,E_\phi~|~i\in\Ob(I),\phi\in\rMorph(I))$ such that :
$$
E_i:=\coprod_{f:i\to t}T^*_fE
\qquad
\text{for every $i\in\Ob(I)$}
$$
where the $f$ ranges over all the morphisms $i\to t$ in $I$.
For every morphism $\phi:i\to j$ of $I$, the morphism
$E_\phi:T^*_\phi E_j\to E_i$ is the unique one fitting in the
commutative diagram :
$$
{\spreaddiagramcolumns{+20pt}\diagram
T^*_\phi T^*_fE \ar[r]^-{\gamma^{T^*}_{\phi,f}}
\ar[d] & T^*_{f\circ\phi}E \ar[d] \\
E_j \ar[r]^-{E_\phi} & E_i
\enddiagram}
\qquad
\text{for every morphism $f:j\to t$ in $I$}
$$
where $\gamma^{T^*}$ is the coherence constraint of $T^*$, and
where the vertical arrows are induced by the chosen universal
cocones for the representatives $E_i$ and $E_j$ of the foregoing
coproducts. For every morphism $h$ in $T_t$, the morphism
$\sTop(T/\rho^t)_!(h)$ is the natural transformation that assigns
to every $i\in\Ob(I)$ the coproduct of morphisms
$\coprod_{f:i\to t}T^*_f(h)$ : details left to the reader.
\end{remark}

\sset\subsubsection{}\label{subsec_total-topos-is-topos}
Let $T$ be a fibred topos over $I$; we consider the associated
fibred lex-site $\pi:\underline\Can(T)\to I$ as in
\eqref{subsec_two-fibrations}, and its total site
$(\underline\Can(T),J)$. Also, for every $t\in\Ob(I)$, we denote
by $i_t:T_t\to\underline\Can(T)$ the inclusion functor.

\begin{theorem}\label{th_top-is-a-topos}
{\em (i)}\ \
With the notation of \eqref{subsec_total-topos-is-topos}, we
have a natural equivalence of categories:
$$
a_T:\sTop(T)\isom(\underline\Can(T),J)^\sim_\sU
$$
and a morphism of sites :
$$
b_T:\Can(\sTop(T))\to(\underline\Can(T),J)
$$
fitting into an essentially commutative diagram (notation
of remark {\em\ref{rem_fibrewise-lims-in-tot-Top}(ii)}) :
$$
{\spreaddiagramcolumns{+10pt}\diagram
\Can(\sTop(T))^\sim_\sU \ar[rd]_{b^\sim_{T*}}
& \sTop(T) \ar[r]^-{\sTop(T/\rho^t)^*} \ar[d]_{a_T}
\ar[l]_-{h_{\sTop(T)}} & T_t \ar[d]^{h_{T_t}} \\
& (\underline\Can(T),J)^\sim_\sU \ar[r]^-{(i_t)^\sim_{\sU*}}
& \Can(T_t)_\sU^\sim
\enddiagram}
\qquad
\text{for every $t\in\Ob(I)$}.
$$

{\em (ii)}\ \ 
The category $\sTop(T)$ is an $\sU$-topos, which we call the
{\em total topos} of the fibred topos $T$.
\end{theorem}
\begin{proof} Let us remark more generally :

\begin{claim}\label{cl_CETA-funckoff}
For every pseudo-functor $\sc:I\to\lex.\Site$ there exists a
natural isomorphism of categories (notation of
\eqref{subsec_T-and-Can} and remarks
\ref{rem_actually-morph-of-sites} and \ref{rem_fibred-site}(v)) :
$$
\alpha_\sc:\sTop(\sT\circ\sc)\isom
\sT(\totSite\circ\underline\cFib_I(\sc)).
$$
\end{claim}
\begin{pfclaim} By proposition \ref{prop_towards-2-yoneda}, we
may assume that $\sc$ is unital, and we let $\gamma^\sc$ be the
coherence constraint of $\sc$. By proposition
\ref{prop_sheaves-on-tot-site}(iv), a presheaf $F$ on the total
site of $\underline\cFib_I(\sc)$ is a sheaf if and only if its
restriction $F_i$ to each fibre site $\sc_i$ is a sheaf. Hence,
such a sheaf $F$ is the datum
$F_\bullet:=(F_i,F_\phi:F_j\to(\sc_\phi)_*^\sim F_i~|~i\in\Ob(I),\ 
(\phi:i\to j)\in\rMorph(I))$ of a system of sheaves and
morphisms of sheaves that make commute the diagrams
$$
\xymatrix{ F_k \ar[d]_-{F_\psi} \ar[rr]^-{F_{\psi\circ\phi}} & &
(\sc_{\psi\circ\phi})_*^\sim F_i \\
(\sc_\psi)_*^\sim F_j \ar[rr]^-{(\sc_\psi)_*^\sim F_\phi} & &
(\sc_\psi)_*^\sim\circ(\sc_\phi)_*^\sim F_i
\ar[u]_{(\gamma^\sc_{\phi,\psi})^\sim_{*,F_i}} 
}$$
for every pair of morphisms
$i\xrightarrow{\phi}j\xrightarrow{\psi}k$ of $I$. Namely,
$F_{\phi,X}:=F(\phi,\one_{\sc_\phi X}):F_jX\to F_i(\sc_\phi X)$
for every such $\phi$, and every $X\in\Ob(T_j)$. Recall
now that, for every site $(\cC,J)$, the topos $\sT(\cC,J)$
is isomorphic to $(\cC,J)^\sim$; especially,
$\sT(\sc_i)$ is isomorphic to the category of sheaves on
the site $\sc_i$, for every $i\in\Ob(I)$. Under this
identification, we then see that a datum $F_\bullet$ as
in the foregoing corresponds precisely to an object of
the category $\sTop(\sT\circ\sc)_*$, and likewise it is
easily seen that the morphisms of $\sTop(\sT\circ\sc)_*$
correspond naturally to the morphisms of sheaves on the
total site of $\underline\cFib_I(\sc)$. Then the sought
isomorphism is the composition of this natural
identification with the isomorphism
$\sTop(\sT\circ\sc)_*\isom\sTop(\sT\circ\sc)_*$
of lemma \ref{lem_change-variance}.
\end{pfclaim}

Recall now that the unit $\eta:\one_{\Topos}\to\sT\circ\Can$ of
the $2$-adjoint pair $(\Can,\sT)$ is a pseudo-natural equivalence 
(theorem \ref{th_adj-topos-site}(iii)); there follows a
pseudo-natural equivalence $\eta_T:T\isom\sT\circ\Can\circ T$.
Then the sought equivalence $a_T$ is the composition of
$\sTop(\eta_T):\sTop(T)\to\sTop(\sT\circ\Can\circ T)$ with the
isomorphism $\alpha_{\Can\circ T}$ of claim \ref{cl_CETA-funckoff}.
The essential commutativity of the square subdiagram of the
diagram of (i) follows by direct inspection.

Now, $(\underline\Can(T),J)$ is isomorphic to an $\sU$-site (remark
\ref{rem_fibred-site}(iii)), hence $(\underline\Can(T),J)^\sim_\sU$
is isomorphic to an $\sU$-topos (remark \ref{rem_summarized}(iv)),
whence (ii). Especially, for every pseudo-functor $\sc:I\to\lex.\Site$
the category $\sTop(\sT\circ\sc)$ is an $\sU$-topos, and then the
$2$-adjunction of theorem \ref{th_adj-topos-site} assigns to
the isomorphism $\alpha_\sc$ of claim \ref{cl_CETA-funckoff}
a morphism of sites :
$$
\beta_\sc:\Can(\sTop(\sT\circ\sc))\to
\totSite\circ\underline\cFib_I(\sc)
$$
so that we may let $b_T:=\beta_{\Can\circ T}\circ\Can(\sTop(\eta_T))$.
The essential commutativity of the triangular subdiagram of (i)
follows from the explicit construction of the unit of this
$2$-adjunction, provided by the proof of theorem \ref{th_adj-topos-site}.
\end{proof}

\begin{remark}\label{rem_describe-a-b}
(i)\ \
By unwinding the constructions in the proof of theorem
\ref{th_top-is-a-topos}, we see that the equivalence $a_T$
assigns to every object
$E:=((E_i,E_\phi)~|~i\in\Ob(I),\phi\in\rMorph(I))$ of
$\sTop(T)$ the datum $(h_{E_i},g_{E_\phi}~|~
i\in\Ob(I),\phi\in\rMorph(I))$, consisting of the sheaf
$h_{E_i}$ on $T_i$ represented by $E_i$ for every $i\in\Ob(I)$,
and the morphism
$$
g_{E_\phi}:h_{E_j}\xrightarrow{h_{E^\dagger_\phi}}h_{T_{\phi_*}E_i}\isom
(T_\phi^*)^\sim_*h_{E_i}
\qquad
\text{for every $(i\xrightarrow{\phi}j)\in\rMorph(I)$}.
$$
Here $E^\dagger_\phi$ is the adjoint to the morphism $E_\phi$,
relative to the adjunction for the pair $(T_\phi^*,T_{\phi*})$
that defines the morphism of topoi $T_\phi:T_i\to T_j$. The
natural identification
$h_{T_{\phi_*}E_i}\isom(T_\phi^*)^\sim_*h_{E_i}$ is deduced as
well from the same adjunction, so the morphism of sheaves
$g_{E_\phi}$ turns out be given by the rule :
$(f:X\to E_j)\mapsto(E_\phi\circ T^*_\phi f)$ for every
$X\in\Ob(S_j)$ and every $f\in h_{E_j}(X)$.

(ii)\ \
On the other hand, the morphism of sites $b_T$ of theorem
\ref{th_top-is-a-topos} is a composition
$$
\underline\Can(T)\xrightarrow{h^a_{\underline\Can(T)}}
(\underline\Can(T),J)^\sim\xrightarrow{a^{-1}_T}\sTop(T)
$$
where we have denoted by $a^{-1}_T$ any choice of a quasi-inverse
of $a_T$. To describe $b_T$ more explicitly, consider any
$(i,X)\in\Ob(\underline\Can(T))$; hence $i\in\Ob(I)$ and
$X\in\Ob(T_i)$; the proof of theorem \ref{th_top-is-a-topos}
assigns to $(i,X)$ the system
$(F(i,X)_j,F(i,X)_\phi~|~j\in\Ob(I),\phi\in\rMorph(I))$
where $F(i,X)_j$ is the restriction to $\Can(T_j)$ of the sheaf
$h^a_{(i,X)}$ on $(\underline\Can(T),J)$. Now, let $G(i,X)_j$ be
the restriction of $h_{(i,X)}$ to $T_j$; for every $Y\in\Ob(T_j)$
we have
$$
G(i,X)_j(Y)=\{(\phi,f)~|~
\phi\in\Hom_I(j,i),f\in\Hom_{T_j}(Y,T^*_\phi X)\}.
$$
For every morphism $h:Y'\to Y$ in $T_j$, the map $G(i,X)_j(h)$
is given by the rule : $(\phi,f)\mapsto(\phi,f\circ h)$. Then set
$$
[i,X,j]:=\coprod_{\phi:j\to i}T^*_\phi X\in\Ob(T_j)
$$
and fix a universal cocone
$(e_{X,\phi}:T^*_\phi X\to[i,X,j]~|~j\xrightarrow{\phi}i)$.
Notice the morphism of presheaves
\set\begin{equation}\label{eq_mansplaining}
G(i,X)_j\to h_{[i,X,j]}
\qquad
(\phi,f:Y\to T^*_\phi X)\mapsto(e_{X,\phi}\circ f:Y\to [i,X,j]).
\end{equation}
Let us check that \eqref{eq_mansplaining} is bicovering.
Indeed, \eqref{eq_mansplaining} is clearly a monomorphism,
hence it suffices to show that it is a covering morphism.
Thus, let $f:Y\to [i,X,j]$ be any morphism on $T_j$ and for
every $\phi\in\Hom_I(j,i)$ set $Y_\phi:=Y\times_{[i,X,j]}T^*_\phi$.
The induced cocone $(e_\phi:Y_\phi\to Y~|~j\xrightarrow{\phi}i)$
is still universal (remark \ref{rem_summarized}(i)) hence it
defines a covering family for $Y$ in the canonical topology
of $T_j$ (remark \ref{rem_summarized}(ii)). But clearly
$f\circ e_\phi$ lies in the image of the map
$G(i,X)_j(Y_\phi)\to h_{[i,X,j]}(Y_\phi)$ for every $\phi:j\to i$,
whence the contention. Hence, \eqref{eq_mansplaining} induces
an isomorphism $G(i,X)_j^a\isom h_{[i,X,j]}$, from the sheaf
$G(i,X)_j^a$ on $\Can(T_j)$ associated to the presheaf $G(i,X)_j$;
the latter is naturally isomorphic to $F(i,X)_j$, according
to proposition \ref{prop_sheaves-on-tot-site}(v). Summing up,
we have exhibited a natural isomorphism of sheaves on $\Can(T_j)$ :
$$
F(i,X)_j\isom h_{[i,X,j]}
\qquad
\text{for every $j\in\Ob(I)$}.
$$

(iii)\ \
Next, according to the proof of theorem \ref{th_top-is-a-topos},
for every morphism $\psi:j\to k$ in $I$, the morphism of sheaves
$F(i,X)_\psi:F(i,X)_k\to(T^*_\psi)^\sim_*F(i,X)_j$ assigns to every
$Y\in\Ob(T_k)$ the map $h^a_{(i,X)}(\psi,\one_{T^*_\psi Y}):h^a_{(i,X)}(k,Y)
\to h^a_{(i,X)}(j,T^*_\psi Y)$. Then, for every such $\phi$ and $Y$ let
$$
G(i,X)_{\psi,Y}:=h_{(i,X)}(\psi,\one_{T^*_\psi Y}):G(i,X)_k(Y)\to
(T^*_\psi)^\wedge G(i,X)_j(Y).
$$
Explicitly, if $\gamma^{T^*}$ denotes the coherence constraint
of the pseudo-functor $T^*$, we get :
$$
G(i,X)_{\psi,Y}(\phi,f)=
(\phi\circ\psi,\gamma^{T^*}_{(\psi,\phi),X}\circ T^*_\phi f)
\qquad
\text{for every $(\phi,f)\in G(i,X)_k(Y)$}.
$$
The rule : $Y\mapsto G(i,X)_{\psi,Y}$ defines a morphism
$G(i,X)_\psi:G(i,X)_k(Y)\to(T^*_\psi)^\wedge G(i,X)_j$ of presheaves,
and $G(i,X)^a_\psi=F(i,X)_\psi$ for every morphism $\psi$ of $I$.
Now, denote by $[i,X,\psi]:T^*_\psi[i,X,k]\to[i,X,j]$ the unique
morphism in $T_j$ fitting into the commutative diagrams :
$$
{\spreaddiagramcolumns{+40pt}\diagram
T^*_\psi T^*_\phi X \ar[r]^-{\gamma^{T^*}_{(\psi,\phi),X}}
\ar[d]_{T^*_\psi(e_{X,\phi})} & T^*_{\phi\circ\psi}X \ar[d]^{e_{X,\phi\circ\psi}} \\
T^*_\psi[i,X,k] \ar[r]^-{[i,X,\psi]} & [i,X,j]
\enddiagram}
\qquad
\text{for every $\phi:k\to i$}.
$$
The foregoing discussion easily implies that we have a
commutative diagram of sheaves :
$$
\xymatrix{ F(i,X)_k \ar[r]^-\sim \ar[d]_{F(i,X)_\psi} &
h_{[i,X,k]} \ar[d]^{h^\dagger_{[i,X,\psi]}} \\
(T^*_\psi)^\sim_*F(i,X)_j \ar[r]^-\sim &
(T^*_\psi)^\sim_*h_{[i,X,j]}
}$$
whose horizontal arrows are the isomorphisms constructed in
(ii), and the adjoint $h^\dagger_{[i,X,\psi]}$ of $h_{[i,X,\psi]}$
is described explicitly as in (i). Summing up, we find that
(up to natural isomorphism) the functor $b_T$ is given by the
rule :
$$
(i,X)\mapsto([i,X,j],[i,X,\psi]~|~j\in\Ob(I),\psi\in\rMorph(I))
\qquad
\text{for every $(i,X)\in\Ob(\underline\Can(T))$}.
$$
\end{remark}

\begin{corollary}\label{cor_excruciating}
Let $I,J$ be two small categories, and  $T,S$ two fibred
topoi over $I$. We have :

{\em (i)}\ \
For every morphism $\omega:T\to S$ of fibred topoi over $I$,
let $\eta^{\sTop(\omega)}$ be the unit of the adjunction for
the pair of functors $(\sTop(\omega)^*,\sTop(\omega)_*)$
provided by remark {\em\ref{rem_top-adjoint-pair}}. Then
$(\sTop(\omega)^*,\sTop(\omega)_*,\eta^{\sTop(\omega)})$
is a morphism of topoi :
$$
\sTop(\omega):\sTop(T)\to\sTop(S).
$$

{\em (ii)}\ \
For every functor $\rho:J\to I$, let $\eta^{\sTop(T/\rho)}$ be
the unit of any adjunction for the pair
$(\sTop(T/\rho)^*,\sTop(T/\rho)_*)$. Then
$(\sTop(T/\rho)^*,\sTop(T/\rho)_*,\eta^{\sTop(T/\rho)})$ is a
morphism of topoi :
$$
\sTop(T/\rho):\sTop(T\circ\rho)\to\sTop(T).
$$

{\em (iii)}\ \
For every pair $\omega:T\to S$, $\nu:S\to U$
of morphisms of fibred topoi over $I$, we have :
$$
\sTop(\nu\circ\omega)=\sTop(\nu)\circ\sTop(\omega).
$$
\end{corollary}
\begin{proof} Assertion (ii) is clear from remark
\ref{rem_fibrewise-lims-in-tot-Top}(i).

(i): In view of proposition \ref{prop_top-adjoint-pair}
and theorem \ref{th_top-is-a-topos}, there remains only to check
that $\sTop(\omega)^*$ is left exact. Arguing as in remark
\ref{rem_use-oldstuff}(i,iii) we may assume that $T$ and
$S$ are unital; then a functor $E_\bullet:J\to\sTop(S)$ assigns
to every $j\in\Ob(J)$ a datum
$E_{\bullet,j}:=(E_{ij},E_{\phi,j}~|~i\in\Ob(I),\ \phi\in\rMorph(I))$
as in remark \ref{rem_use-oldstuff}(i), and to every morphism
$f:j\to j'$ a system of morphisms
$(E_{i,f}:E_{ij}\to E_{ij'}~|~i\in\Ob(I))$ such that
$$
E_{\phi,j'}\circ S^*_\phi E_{i',f}=E_{i,f}\circ E_{\phi,j}
\qquad
\text{for every morphism $\phi:i\to i'$ of $I$}.
$$
The limit of $E_\bullet$ in the category $\sTop(S)$ is
the same as the colimit of $E^o_\bullet$ in the category
$\Sigma(\cFib(S^*)/I)$, and we construct the latter following
remark \ref{rem_fibrewise-lims-in-tot-Top}(iii). Thus,
for every $i\in\Ob(I)$ we consider the functor
$(E^o_\bullet)^i:J\to\cFib(S^*)$ given by the rules :
$j\mapsto E_{ij}$ for every $j\in\Ob(J)$ and
$f\mapsto E^o_{i,f}$ for every morphism $f$ of $J$. Then
$(E^o_\bullet)^i$ factors through the inclusion functor
$S^o_i\to\cFib(S^*)$, and the colimit of the resulting
functor $(\bar E{}^o_\bullet)^i:J\to S^o_i$ in $S^o_i$ is
the limit of the opposite functor
$\bar E{}_\bullet^i:J^o\to S_i$, so we pick
$L(i)\in\Ob(S_i)$ representing the latter limit and a
universal cone $\tau^i:c_{L(i)}\Rightarrow E_\bullet^i$.
Next, every morphism $\phi:i\to i'$ in $I$ induces a
natural transformation
$$
\bar E{}^\phi_\bullet:
S^*_\phi\circ\bar E{}^{i'}_\bullet\Rightarrow\bar E{}^i_\bullet
\qquad
j\mapsto E_{\phi,j}:S^*_\phi E_{i'j}\to E_{ij}.
$$
Then there exists a unique morphism
$L(\phi):S^*_\phi L(i')\to L(i)$ of $S_i$ fitting into the
commutative diagram :
$$
\cD
\qquad : \qquad
{\diagram S^*_\phi\circ c_{L(i')} \ar@{=>}[rr]^-{c_{L(\phi)}}
\ar@{=>}[d]_{S^*_\phi*\tau^{i'}} & & c_{L(i)} \ar@{=>}[d]^{\tau^i} \\
S^*_\phi\circ\bar E{}^{i'}_\bullet \ar@{=>}[rr]^-{\bar E{}^\phi_\bullet} & &
\bar E{}^i_\bullet
\enddiagram}$$
and remark \ref{rem_fibrewise-lims-in-tot-Top}(iii) shows that
$L_\bullet:=(L(i),L(\phi)~|~i\in\Ob(I)\, \phi\in\rMorph(I))$ is
an object of $\sTop(S)$ representing the limit of $E_\bullet$,
and the system of cones $(\tau^i~|~i\in\Ob(I))$ adds up to
a universal cone $\tau:c_{L_\bullet}\Rightarrow E_\bullet$.
Now, set
$$
F_\bullet:=\sTop(\omega)^*\circ E_\bullet:J\to\sTop(T).
$$
We define likewise the functors $\bar F{}^i_{\!\!\bullet}:J\to T_i$
for every $i\in\Ob(I)$ and the natural transformations
$\bar F{}^\phi_{\!\!\bullet}:
T^*_\phi\circ\bar F{}^{i'}_{\!\!\bullet}\Rightarrow\bar F{}^i_{\!\!\bullet}$
for every morphism $\phi:i\to i'$ of $I$; then by a simple
inspection we deduce from $\cD$ the commutative diagram
$$
\cD^\omega
\qquad : \qquad
{\diagram T^*_\phi\circ c_{\omega^*_{i'}L(i')}
\ar@{=>}[rr]^-{c_{L(\phi)^\omega}}
\ar@{=>}[d]_{T^*_\phi*\omega^*_{i'}*\tau^{i'}} & & c_{\omega^*_iL(i)}
\ar@{=>}[d]^{\omega^*_i*\tau^i} \\
T^*_\phi\circ\bar F{}^{i'}_{\!\!\bullet}
\ar@{=>}[rr]^-{\bar F{}^\phi_{\!\!\bullet}} & &
\bar F{}^i_{\!\!\bullet}
\enddiagram}$$
with
$(\omega^*_iL(i),L(\phi)^\omega~|~i\in\Ob(I),\ \phi\in\rMorph(I)):=
\sTop(\omega)^*L_\bullet$ as in remark \ref{rem_use-oldstuff}(ii).
Lastly, if $J$ is a finite category, then by assumption
$\omega^*_i*\tau^i$ is still a universal cone for every $i\in\Ob(I)$,
and therefore the commutativity of $\cD^\omega$ determines
again $L(\phi)^\omega$ uniquely. In this case, remark
\ref{rem_fibrewise-lims-in-tot-Top}(iii) shows that
$\sTop(\omega)^*L_\bullet$ represents the limit of $F_\bullet$
and $\sTop(\omega)^**\tau$ is a universal cone.

(iii) follows easily from the explicit description of
$\eta^{\sTop(\omega)},\eta^{\sTop(\nu)}$ and $\eta^{\sTop(\nu\circ\omega)}$
given by remark \ref{rem_top-adjoint-pair}.
\end{proof}

\begin{remark}\label{rem_b-bullet}
Let $I$ be any small category.

(i)\ \
From corollary \ref{cor_excruciating}(i,iii) it follows
easily that the rules : $T\mapsto\sTop(T)$,
$\omega\mapsto\sTop(\omega)$ and $\Xi\mapsto\sTop(\Xi)_*$
for every fibred topos $T$ over $I$, every morphism $\omega$
of such fibred topoi, and every transformation $\Xi$ of
such morphisms, define a strict pseudo-functor
$$
\sTop:\sPsFun(I,\Topos)\to\Topos.
$$

(ii)\ \
It is also easily seen that the isomorphisms of claim
\ref{cl_CETA-funckoff} yield a pseudo-natural isomorphism
of pseudo-functors :
$$
\alpha_\bullet:\sTop\circ\sPsFun(I,\lex.\sT)\isom
\sT\circ\totSite\circ\lex.\underline\cFib_I.
$$
Moreover, let
$\underline\sc^\bullet:\fib.\lex.\Site(I)\to\sPsFun(I,\lex.\Site)$
be a strict and strong pseudo-inverse for the $2$-equivalence
$\lex.\underline\cFib_I$ (see remark \ref{rem_fibred-site}(v));
then from $\alpha_\bullet*\underline\sc^\bullet$ we deduce a
pseudo-natural isomorphism of pseudo-functors (notation of
\eqref{subsec_two-fibrations}) :
$$
\sa_\bullet:\sTop\circ\underline{\lex.\sT}\isom\sT\circ\totSite.
$$
By $2$-adjunction, from $\alpha_\bullet$ and $\sa_\bullet$
we then deduce as well pseudo-natural transformations
$$
\begin{aligned}
\beta_\bullet&\,:\Can\circ\sTop\circ\sPsFun(I,\lex.\sT)\to
\totSite\circ\lex.\underline\cFib_I\qquad\sc\mapsto\beta_\sc \\
\sbb_\bullet&\,:\Can\circ\sTop\circ\underline{\lex.\sT}\to\totSite
\end{aligned}
$$
where $\beta_\sc$ is defined as in the proof of theorem
\ref{th_top-is-a-topos} : the details are left to the reader.

(iii)\ \
Lastly, the rule $T\mapsto a_T$ of theorem
\ref{th_top-is-a-topos}(i) defines the pseudo-natural equivalence :
$$
a_\bullet:=\sPsFun(I,\eta)\odot(\sa_\bullet*\underline\Can):
\sTop\isom\sT\circ\totSite\circ\underline\Can
$$
where $\eta:\one_{\Topos}\to\sT\circ\Can$ is the unit of
the $2$-adjoint pair $(\Can,\sT)$. Again by $2$-adjunction,
the rule : $T\mapsto b_T$ yields likewise a pseudo-natural
transformation
$$
b_\bullet:\Can\circ\sTop\to\totSite\circ\underline\Can.
$$
\end{remark}

\subsection{Localization and points of a topos}
\label{sec_Localization-topoi}
Let $\cC$ be a category, and $f:Y\to X$ any morphism of $\cC$.
We consider the source functor $\ss_X:\cC\!/\!X\to\cC$ as in
\eqref{subsec_slice-cat}, and the induced functor
$$
(\ss_X)_{|f}:(\cC/X)/f\to\cC/(\ss_X f)=\cC/Y
$$
as in \eqref{eq_restrict-over-X}, which is obviously an
isomorphism of categories, hence it induces a bijection :
$$
\{\text{sieves of $\cC/Y$}\}\isom
\{\text{sieves of $(\cC/X)/f$}\}
\qquad
\cS\mapsto(\ss_X)^{-1}_{|f}\cS
$$
(notation of definition \ref{def_sieve}(iii)).
Fix a universe $\sV$ with $\sU\subset\sV$ and such that
$\cC$ is $\sV$-small, so that the functor
$(\ss_X)_{\sV!}:(\cC/X)_\sV^\wedge\to\cC_\sV^\wedge$ is well
defined. We notice that for every subobject $R$ of the
presheaf $h_f$ on $\cC/X$, the presheaf $(\ss_X)_{\sV!}R$
is a subobject of $h_Y$ : more precisely, in light of
\eqref{eq_can-be-used} we see that for a sieve $\cS$
of $(\cC/X)/f$ and a sieve $\cT$ of $\cC/Y$, we have :
\set\begin{equation}\label{eq_equivalence-sieve}
h_\cT=(\ss_X)_{\sV!}h_\cS
\quad\Leftrightarrow\quad
\cS=(\ss_X)^{-1}_{|f}\cT.
\end{equation}
Let now $J$ be a topology on $\cC$; endow $\cC/X$ with the
topology $J_X$ induced by $J$ via $\ss_X$, and set 
$$
C:=(\cC,J)
\qquad
C/X:=(\cC/X,J_X).
$$
Since $(\ss_X)_{\sV!}$ commutes with fibre products (proposition
\ref{prop_in-the-same-vein}(vi.c)), lemma
\ref{lem_induced-top}(ii.a) says that a subobject $R\subset h_f$
covers $f$ for the topology $J_X$ if and only if the induced morphism
$(\ss_X)_{\sV!}R\to(\ss_X)_{\sV!}(h_f)=h_Y$ covers $Y$ for the topology
$J$. In other words, a family of morphisms
$(h_\lambda/X:(f_\lambda:Y_\lambda\to Y)\to f~|~\lambda\in\Lambda)$
covers $f$ for the topology $J_X$ if and only if the family
$(h_\lambda:Y_\lambda\to Y~|~\lambda\in\Lambda)$ covers $Y$
for the topology $J$. In view of \eqref{eq_equivalence-sieve},
it follows easily that $\ss_X$ is both continuous and cocontinuous
for the topologies $J$ and $J_X$.

\sset\subsubsection{}\label{subsec_topol-on-C-over-X}
Suppose now that $C$ is a $\sU$-site; then the same holds
for $C/X$ as well : indeed, if $G\subset\Ob(\cC)$ is a small
topologically generating family for $C$, then $G/X\subset\Ob(\cC/X)$
is a small topologically generating family for $C/X$
(notation of definition \ref{def_U-site}(i)). By virtue of
corollary \ref{cor_two-U-sites}, the functor
$\tilde\ss_{X*}\simeq\breve\ss{}_X^*:C^\sim\to(C/X)^\sim$ admits
therefore both a left adjoint $\tilde\ss{}_X^*$ and a right
adjoint $\breve\ss_{X*}$. We introduce a special notation and
terminology for these functors :
\begin{itemize}
\item
The functor $\tilde\ss_{X*}$ shall be also denoted
$j_X^*$, and called the functor of {\em restriction to $X$}.
\item
The functor $\breve\ss_{X*}$ shall be denoted $j_{X*}$, and
called the {\em direct image\/} functor.
\item
The functor $\tilde\ss{}_X^*$ shall be denoted $j_{X!}$,
and called the functor of {\em extension by empty}.
\end{itemize}
Thus, $j^*_X$ is right adjoint to $j_{X!}$, and left adjoint to
$j_{X*}$. However, a left adjoint for $j^*_X$ can be exhibited
alternatively by a more explicit construction, which will allow
us to extract some useful additional properties of the category
$(C/X)^\sim$. Indeed, even though $\cC/X$ is not necessarily
$\sU$-small, we can invoke proposition \ref{prop_in-the-same-vein}(vi.a)
and remark \ref{rem_U-site}(ii), in order to define a functor
by the same rule as in \eqref{subsec_both-pigeons}, namely
$$
j_{X!}:(C/X)_\sU^\sim\to C_\sU^\sim
\qquad
F\mapsto(\ss_{X!}\circ i_{C/X}F)^a
$$
where $i_{C/X}:(C/X)^\sim\to(\cC/X)^\wedge$ is the forgetful functor.
Moreover, remark \ref{rem_rep-and-sheafify}(v) (and again
remark \ref{rem_U-site}(ii)) implies that this functor is
(isomorphic to) the restriction of $(\tilde\ss_X)^*_\sV$, and
since the inclusion functor $C^\sim_\sU\to C^\sim_\sV$ is fully faithful,
we deduce that it is also a left adjoint to $j^*_X$. Furthermore,
recall that $\ss_{X!}$ factors through the source functor
$\ss_{h_X}:\cC^\wedge/h_X\to\cC^\wedge$ and an equivalence
$e_X:(\cC/X)^\wedge\isom\cC^\wedge/h_X$ (proposition
\ref{prop_in-the-same-vein}(vi.b)); consequently, $j_{X!}$
factors through the source functor $\ss_{h_X^a}:C^\sim/h_X^a\to C^\sim$
and the composition
$$
\tilde e_X:(C/X)^\sim\xrightarrow{i_{C/X}}(\cC/X)^\wedge
\xrightarrow{e_X}\cC^\wedge/h_X\xrightarrow{(-)^a_{|h_X}}C^\sim\!/h_X^a
$$
where $(-)^a_{|h_X}$ is induced by $(-)^a:\cC^\wedge\to C^\sim$ and
the object $h_X\in\Ob(\cC^\wedge)$, as in \eqref{eq_restrict-over-X}.

\begin{remark}\label{rem_continue-local}
(i)\ \
Let $g:Y\to Z$ be any morphism in $\cC$; then we have the
sites $C/Y$ and $C/Z$ as in \eqref{sec_Localization-topoi},
as well as the functor $g_*:\cC/Y\to\cC/Z$ of \eqref{eq_push-for}.
The isomorphism of categories $(\ss_Z)_{|g}:(\cC/Z)/g\isom\cC/Y$
identifies $g_*$ with the source functor
$\ss_g:(\cC/Z)/g\to\cC/Z$. Hence, the discussion of
\eqref{subsec_topol-on-C-over-X} applies to $g_*$ as well,
and since $\ss_Z\circ g_*=\ss_Y$, we see that $J_Y$ is the
topology induced by $J_Z$ via $g_*$ on $\cC/Y$; so, $g_*$
is continuous and cocontinuous for the topologies $J_Y$
and $J_Z$.  Moreover, if $C$ is a $\sU$-site, then
$(\tilde g_*)_*:(C/Z)^\sim\to(C/Y)^\sim$ admits a left adjoint
$(\tilde g_*)^*$ and a right adjoint $(\breve g_*)_*$. In
agreement with the foregoing, we let
$$
j_g^*:=(\tilde g_*)_*
\qquad
j_{g*}:=(\breve g_*)_*
\qquad
j_{g!}:=(\tilde g_*)^*.
$$
Clearly $j_g^*\circ j_Z^*=j_Y^*$, and we have isomorphisms of
functors : $j_{Z*}\circ j_{g*}\isom j_{Y*}$ and
$j_{Z!}\circ j_{g!}\isom j_{Y!}$.

(ii)\ \
By propositions \ref{prop_half-dominates}(i) and
\ref{prop_was-get-maddd}(iii,iv), we also see that the pairs
$(j^*_X,j_{X*})$ and $(j_g^*,j_{g*})$ determine morphisms of
topoi, unique up to unique isomorphism :
$$
j_X:(C/X)^\sim\to C^\sim
\qquad
j_g:(C/Y)^\sim\to(C/Z)^\sim.
$$

(iii)\ \
Suppose that all finite products of $\cC$ are representable.
Then for every $X\in\Ob(\cC)$, the source functor $\ss_X$ admits
a right adjoint
$$
p_X:\cC\to\cC/X
\qquad
Y\mapsto(p_{X,Y}:X\times Y\to X)
$$
where $X\times Y$ is any choice of a representative for the
product of $X$ and $Y$, and $p_{X,Y}$ is the corresponding natural
projection : see proposition \ref{prop_in-the-same-vein}(iii).
Since $\ss_X$ is cocontinuous for the topologies $J$ and $J_X$,
the functor $p_X$ is a morphism of sites $C/X\to C$ (remark
\ref{rem_choose-two-univs}(v)), and we have an isomorphism
of morphisms of topoi :
$$
j_X\isom\tilde p_X.
$$
\end{remark}

\begin{proposition}\label{prop_localize-topos}
With the notation of \eqref{subsec_topol-on-C-over-X}, the
following holds :
\begin{enumerate}
\item
The functor $\tilde e_X$ is an equivalence.
\item
We have essentially commutative diagrams :
$$
\xymatrix@C+40pt{ (C/X)^\sim \ar[r]^-{i_{C/X}} \ar[d]_{\tilde e_X} &
(\cC/X)^\wedge \ar[d]^{e_X} &
(\cC/X)^\wedge \ar[d]_{e_X} \ar[r]^-{(-)^a} &
(C/X)^\sim  \ar[d]^{\tilde e_X} \\
C^\sim/h^a_X \ar[r]^-{(i_C)_{|h_X}} & \cC^\wedge/h_X &
\cC^\wedge/h_X \ar[r]^-{(-)^a_{|h_X}} & C^\sim/h^a_X
}$$
where $(i_C)_{|h_X}$ and $(-)^a_{|h_X}$ are the functors attached
to the adjoint pair $((-)^a,i_C)$ and the object
$h_X\in\Ob(\cC^\wedge)$, as in example {\em\ref{ex_comma-adjunction}(i)}.
\end{enumerate}
\end{proposition}
\begin{proof}Let $\eta_F:h_X\to h^a_X$ be the natural
morphism in $\cC^\wedge$. We remark :

\begin{claim}\label{cl_cazzuto}
For every presheaf $F$ on $\cC/X$ there exists a cartesian diagram
in $\cC^\wedge$ :
$$
{\diagram
\ss_{X!}(F^a) \ar[r] \ar[d]_{e_X(F^a)} & (\ss_{X!}F)^a
\ar[d]^{e_X(F)^a} \\
h_X \ar[r]^-{\eta_X} & h^a_X.
\enddiagram}
$$
\end{claim}
\begin{pfclaim} As already pointed out in the proof of
proposition \ref{prop_in-the-same-vein}(vi.b), the presheaf
 $\bone_{\cC/X}:=h_{\one_X}$ is a final object of $(\cC/X)^\wedge$,
and the equivalence $e_X$ assigns to every presheaf $F$ on
$\cC/X$ the $h_X$-presheaf
$e_X(F):\ss_{X!}F\to\ss_{X!}(\bone_{\cC/X})\isom h_X$
deduced from the unique morphism $F\to\bone_{\cC/X}$.
Moreover, $e_X$ admits the quasi-inverse
$$
e'_X:\cC^\wedge/h_X\to(\cC/X)^\wedge
\qquad
(\phi:G\to h_X)\mapsto\ss^\wedge_X(G)\times_{\ss^\wedge_X(h_X)}\bone_{\cC/X}
$$
where the fibre product is defined via the unit of adjunction
$\eta_\bone:\bone_{\cC/X}\to\ss^\wedge_X(\ss_{X!}
\bone_{\cC/X})\isom\ss^\wedge_X(h_X)$ and the morphism
$\ss^\wedge_X(\phi)$. Thus, we have a natural isomorphism
$$
F\isom\ss^\wedge_X(\ss_{X!}F)\times_{\ss^\wedge_X(h_X)}\bone_{\cC/X}
\qquad
\text{in $(\cC/X)^\wedge$}.
$$
On the other hand, the functor $\ss^\wedge_X$ commutes with the
functor $G\mapsto G^a$ (corollary \ref{cor_two-U-sites}(iii)),
so there follows an isomorphism
$$
F^a\isom\ss^\wedge_X(\ss_{X!}F)^a\times_{\ss^\wedge_X(h^a_X)}\bone_{\cC/X}
\qquad
\text{in $(C/X)^\sim$}.
$$
But it is easily seen that the morphism
$\eta^a_\bone:\bone_{\cC/X}\to\ss^\wedge_X(h^a_X)$ is the composition
of $\eta_\bone$ and the natural morphism
$\ss^\wedge_X(h_X)\to\ss^\wedge_X(h^a_X)$, so we get a commutative
diagram with cartesian squares:
$$
\xymatrix{ F^a \ar[r]^-\alpha \ar[d] &
\ss^\wedge_X((\ss_{X!}F)^a\times_{h^a_X}h_X) \ar[r] \ar[d] &
\ss^\wedge_X(\ss_{X!}F)^a \ar[d] \\
\bone_{\cC/X} \ar[r] & \ss^\wedge_X(h_X) \ar[r] & \ss^\wedge_X(h^a_X).
}$$
Especially, $\alpha$ induces an isomorphism
$F^a\isom e'_X((\ss_{X!}F)^a\times_{h^a_X}h_X)$, whence an
isomorphism
$$
e_X(F^a)\isom((\ss_{X!}F)^a\times_{h^a_X}h_X\to h_X)
$$
and the claim follows.
\end{pfclaim}

(i): Let $\phi:G\to h^a_X$ be any object of $C^\sim/h^a_X$;
set $G':=G\times_{h^a_X}h_X$, and let $\phi':G'\to h_X$ be
the induced projection. Let also $F:=e'_X(\phi')$, which
is a presheaf on $\cC/X$. We claim that $F$ is a sheaf for
the topology $J_X$. Indeed, by claim \ref{cl_cazzuto} we
have an isomorphism in $\cC^\wedge/h_X$ :
$$
\ss_{X!}(F^a)\isom(\ss_{X!}F)^a\times_{h^a_X}h_X\isom
G'^a\times_{h^a_X}h_X\isom G'
$$
and since $e'_X$ is quasi-inverse to $e_X$, we have an
isomorphism of $h_X$-presheaves $\ss_{X!}(F)\isom G'$.
Thus, $F\isom F^a$ in $(\cC/X)^\wedge$, as required. Next
we claim that the resulting functor
$$
\tilde e'_X:C^\sim/h^a_X\to(C/X)^\sim
\qquad
(\phi:G\xrightarrow{\phi}h^a_X)\mapsto
e'_X(G\times_{h^a_X}h_X\xrightarrow{\phi\times_{h^a_X}h_X}h_X)
$$
is a quasi-inverse for $\tilde e_X$. Indeed, for every
$\phi$ as in the foregoing we have an isomorphism
$$
\tilde e_X\circ\tilde e'_X(\phi)\isom(\phi\times_{h^a_X}h_X)^a
\isom\phi
\qquad
\text{in $C^\sim/h^a_X$}
$$
which is natural with respect to $\phi$. Lastly, for every
sheaf $F$ on $C/X$ we have $\tilde e'_X\circ\tilde e_X(F)=
e'_X((e_X(F)^a\times_{h^a_X}h_X)$, which is isomorphic to
$e'_X(e_X(F))\simeq F$, by claim \ref{cl_cazzuto}.

(ii): For every sheaf $F$ on $C/X$, claim \ref{cl_cazzuto}
yields an isomorphism :
$$
e_X(i_{C/X}F)\isom(\tilde e_X(F)\times_{h^a_X}h_X\to h_X).
$$
The essential commutativity of the left square diagram in (ii)
is an immediate consequence; since $e_X$ and $\tilde e_X$ are
both equivalences, and since $(-)^a_{|h_X}$ is left adjoint to
$(i_C)_{|h_X}$, we deduce the essential commutativity also for
the right square diagram.
\end{proof}

\begin{remark}
(i)\ \
Under the equivalence $\tilde e_X$ of proposition
\ref{prop_localize-topos}, the functor $j_X^*$ is identified
with the functor :
$$
C^\sim\to C^\sim/h_X^a
\qquad
F\mapsto(F\times h^a_X\to h^a_X)
$$
and $j_{X!}$ is identified with the source functor
$\ss_{h^a_X}:C^\sim/h^a_X\to C^\sim$ of \eqref{subsec_slice-cat}.

(ii)\ \
Likewise, the functor $j_g^*$ of remark \ref{rem_continue-local}(i)
is identified with the functor
$$
C^\sim/h^a_Z\to C^\sim/h^a_Y
\qquad
(F\to h^a_Z)\mapsto(F\times_{h^a_Z}h^a_Y\mapsto h^a_Y)
$$
and $j_{g!}$ is identified with the functor
$(h_g^a)_*:C^\sim/h^a_Y\to C^\sim/h^a_Z$ induced by
$h^a_g:h^a_Y\to h^a_Z$.
\end{remark}

\begin{proposition}\label{prop_localize-continuity}
Let $C:=(\cC,J)$ and $C':=(\cC';J')$ be two sites, $u:\cC\to\cC'$
a continuous functor, $X$ any object of\/ $\cC$. Then we have :
\begin{enumerate}
\item
The functor $u_{|X}:\cC/X\to\cC'/uX$ is continuous
for the sites $C/X$ and $C'/uX$.
\item
We have essentially commutative diagrams :
$$
\xymatrix@C+20pt{
(C'/uX)^\sim \ar[r]^-{\tilde{u_{|X}}{}_*} \ar[d]_{\tilde e_{uX}} &
(C/X)^\sim \ar[d]^{\tilde e_X} &
(C/X)^\sim \ar[d]_{\tilde e_X} \ar[r]^-{\tilde{u_{|X}}{}^*} &
(C'/uX)^\sim \ar[d]^{\tilde e_{uX}} \\
C'^\sim/h^a_{uX} \ar[r]^-{(\tilde u_*)_{|h^a_X}} &
C^\sim/h^a_X & C^\sim/h^a_X \ar[r]^-{(\tilde u^*)_{|h^a_X}} &
C'^\sim/h^a_{uX}.
}$$
\item
If $u$ is a morphism of sites $C'\to C$, then $u_{|X}$ is a
morphism of sites $C'/uX\to C/X$.
\end{enumerate}
\end{proposition}
\begin{proof}(i): Combining propositions \ref{prop_localize-topos}(ii)
and \ref{prop_done-right}(i) with example \ref{ex_comma-adjunction}(iii)
we get an essentially commutative diagram, whose vertical arrows
are equivalences :
$$
\xymatrix@C+40pt{ (C'/uX)^\sim \ar[r]^-{(u_{|X})^\wedge\circ i_{C'/X}}
\ar[d]_{\tilde e_{uX}} & (\cC/X)^\wedge \ar[d]^{e_X} \\
C'^\sim/h^a_{uX} \ar[r]^-{(u^\wedge\circ i_{C'})_{|h_X}} & \cC^\wedge/h_X.
}$$
Now, let $F$ be a sheaf on $C'/uX$, and denote by $\phi:G\to h^a_{uX}$
the $h^a_{uX}$-sheaf $\tilde e_{uX}(F)$. By proposition
\ref{prop_localize-topos}(ii), it suffices to show that
$(u^\wedge\circ i_{C'})_{|h_X}(\phi)\simeq
(u^\wedge G\times_{u^\wedge h^a_{uX}}h_X\to h_X)$ lies in the essential
image of $(i_C)_{|h_X}$. But the morphism $h_X\to u^\wedge h^a_{uX}$
is the composition of the natural morphisms $h_X\to h^a_X$ and
$h^a_X\to u^\wedge h^a_{uX}$, hence
$u^\wedge G\times_{u^\wedge h^a_{uX}}h_X\isom G'\times_{h^a_X}h_X$, with
$G':=u^\wedge G\times_{u^\wedge h^a_{uX}}h^a_X$, and the latter is a
sheaf on $C$, since $u$ is continuous, whence the contention.

(ii): Recall that we have a natural isomorphism
$\tilde u{}^*(h^a_X)\isom h^a_{uX}$ in $C'^\sim$ (lemma
\ref{lem_cont-funct-site}(ii)). By example
\ref{ex_comma-adjunction}(i), we deduce an adjoint pair
of functors :
$$
(\tilde u{}^*)_{|h^a_X}:C^\sim/h^a_X\to C'^\sim/h^a_{uX}
\qquad
(\tilde u_*)_{|h^a_X}:C'^\sim/h^a_{uX}\to C^\sim/h^a_X.
$$
Now, the proof of (i) yields a natural isomorphism of functors :
$$
\begin{aligned}
\tilde e_X\circ(u_{|X})^\sim_*&\,\isom
(-)^a_{|h_X}\circ e_X\circ i_{C/X}\circ(u_{|X})^\sim_* \\
&\,\isom(-)^a_{|h_X}\circ e_X\circ(u_{|X})^\wedge\circ i_{C'/uX} \\
&\,\isom(-)^a_{|h_X}\circ(u^\wedge\circ i_{C'})_{|h_X}\circ\tilde e_{uX}
\end{aligned}
$$
and a direct inspection of the constructions shows that
$(-)^a_{|h_X}\circ(u^\wedge\circ i_{C'})_{|h_X}$ is isomorphic
to the functor $(\tilde u_*)_{|h^a_X}$, so that the left
square diagram of (ii) is essentially commutative. Since
$\tilde e_X$ and $\tilde e_{uX}$ are equivalences, we deduce
that the same holds also for the right square diagram.

(iii): By assumption, the functor $\tilde u{}^*$ is left
exact, so the same holds for $(\tilde u{}^*)_{|h^a_X}$
(proposition \ref{prop_done-right}(ii)), and by (ii) the
same follows for $\tilde{u_{|X}}{}^*$, whence the contention.
\end{proof}

\begin{example}\label{ex_localization-topos}
(i)\ \
Let $T$ be a topos, and $U\in\Ob(T)$ any object. By proposition
\ref{prop_localize-topos}(i) and remark
\ref{rem_continue-local}(iii), we have a natural equivalence
$$
(\Can(T)/U)^\sim\isom\Can(T)^\sim/h_U
$$
identifying the functor $j^*_U:\Can(T)^\sim\to(\Can(T)/U)^\sim$
with $p_{h_U}:\Can(T)^\sim\to\Can(T)^\sim/h_U$. The source functor
$\ss_{h_U}:\Can(T)^\sim/h_U\to\Can(T)^\sim$ is a {\em left} adjoint
for $p_{h_U}$, but the latter admits also a {\em right} adjoint, which
is naturally identified with $j_{U*}:(\Can(T)/U)^\sim\to\Can(T)^\sim$.
On the other hand, the Yoneda equivalence $T\isom\Can(T)^\sim$
induces an equivalence $T/U\isom\Can(T)^\sim/h_U$, which shows
that $T/U$ is a topos, and identifies $p_{h_U}$ in turn with the
functor
$$
p_U:T\to T/U
\qquad
X\mapsto X_{|U}:=(X\times U\to U).
$$
Hence, the latter is a morphisms of $\sU$-sites
$\Can(T)\to\Can(T/U)$, and we denote again by
$$
j_U:=(j_U^*,j_{U*},\eta^U):T/U\to T
$$
the corresponding morphism of topoi. Moreover, the source
functor $\ss_U$ is a left adjoint for $j^*_U$, and shall
be denoted also by $j_{U!}:T/U\to T$. Likewise, any morphism
$f:U'\to U$ in $T$ determines, up to unique isomorphism, a
morphism of topoi
$$
j_f:=(j_f^*,j_{f*},\eta^f):T/U'\to T/U
$$
with an isomorphism $j_U\circ j_f\isom j_{U'}$ of morphisms of topoi.

(ii)\ \
Recall that the target functor $\st:\sMorph(T)\to T$ is a
fibration (example \ref{ex_fibred-cats}(iii)); pick a
cleavage for $\st$, and let $\sc:T^o\to\sV\tdu\bCat$
be its associated pseudo-functor (for a universe $\sV$
such that $T$ is $\sV$-small). We may choose the cleavage
so that for every morphism $f:U'\to U$, the functor
$\sc_f:T/U\to T/U'$ agrees with $j^*_f$ (the details are
left to the reader). Notice then that
$\sc^o:T\to\sV\tdu\bCat^o$ factors through the forgetful
strict pseudo-functor
$$
\sLink(\sV\tdu\bCat^o)\to\sV\tdu\bCat^o
\qquad
\cA\mapsto\cA
\qquad
(F,G,\eta,\eps)\mapsto F
$$
and a pseudo-functor
$$
T\to\sLink(\sV\tdu\bCat^o)
\qquad
U\mapsto T/U
\qquad
(U'\xrightarrow{f}U)\mapsto(j^*_f,j_{f*},\eta^f,\eps^f)
$$
where $\eps^f$ is the counit for the adjunction of the
pair $(j^*_f,j_{f*})$, whose unit is $\eta^f$. We compose
the latter with the strict isomorphism of $2$-categories
$\sLink(\sV\tdu\bCat^o)\isom{}^o\sLink(\sV\tdu\bCat)$
provided by proposition \ref{prop_isoms-of-links}, and
notice that the resulting pseudo-functor maps $T$ to the
strong sub-$2$-category ${}^o\Topos$ of ${}^o\sLink(\sV\tdu\bCat)$.
Thus, we obtain a well defined fibred topos over $T$ :
$$
T/-:T={}^oT\to\Topos
\qquad
U\mapsto T/U
\qquad
(U'\xrightarrow{f}U)\mapsto j_f.
$$

(iii)\ \
In case $U$ is a subobject of the final object $1_T$ of $T$,
the morphism $j_U$ of (i) is called an {\em open subtopos\/}
of $T$. In this case we denote by $\sC U$ the full subcategory
of $T$ such that
$$
\Ob(\sC U)=\{X\in\Ob(T)~|~j^*_UX=1_{T/U}\}.
$$
Then $\sC U$ is a topos, called the {\em complement of $U$
in $T$}, and the inclusion functor $i_*:\sC U\to T$ admits
a left adjoint $i^*:T\to\sC U$, namely, the functor which
assigns to every $X\in\Ob(T)$ the push-out $X_{|\sC U}$ in
the cocartesian diagram :
$$
\xymatrix{ X\times U \ar[r]^-{p_X} \ar[d]_{p_U} & X \ar[d] \\
                   U \ar[r] & X_{|\sC U}
}$$
where $p_X$ and $p_U$ are the natural projections. Moreover,
$i^*$ is an exact functor, hence the adjoint pair $(i^*,i_*)$
defines a morphism of topoi $\sC U\to T$, unique up to unique
isomorphism. (See \cite[Exp.IV, Prop.9.3.4]{SGA4-1}.)
\end{example}

\begin{remark}\label{rem-sections-of-f}
(i)\ \
Let $f:T'\to T$ be a morphism of topoi, $U\in\Ob(T)$, and let
us fix final objects $1_T$ of $T$ and $1_{T'}$ of $T'$. For
every $Y\in\Ob(T')$, let also $u_Y:Y\to 1_{T'}$ be the unique
morphism in $T'$. Suppose that $s:T'\to T\!/U$ is a morphism
of topoi with an isomorphism $\omega:f\isom j_U\circ s$. Hence,
we have a functorial isomorphism (notation of remark
\ref{rem_adjoint-transf}(ii))
$$
\omega^\dagger_X:s^*(j_U^*X)\isom f^*X
\qquad
\text{for every $X\in\Ob(T)$}.
$$
So, $s^*\one_U\isom f^*1_T$ is a final object of $T'$, hence
$u_{s^*\one_U}:s^*\one_U\to 1_{T'}$ is an isomorphism.
Let $\Delta_U:\one_U\to(U\times U\xrightarrow{j^*_U(U)}U)$ be
the diagonal morphism; then
$$
\sigma(s,\omega):=
\omega^\dagger_U\circ s^*(\Delta_U)\circ u^{-1}_{s^*\one_U}:
1_{T'}\to f^*U
$$
is an element of $\Gamma(T',f^*U)$ (notation
\eqref{subsec_glob-sections}). Moreover, for every object
$\phi:X\to U$ of $T\!/U$, we have a cartesian diagram in $T\!/U$ :
$$
\cD
\quad :\quad
{\diagram
X \ar[rr]^-{\Gamma_{\!\phi}} \ar[dd]_{\phi/U} \ar[rd]^-\phi & &
X\times U \ar[dd]^{j^*_U\phi} \ar[ld]_-{j^*_UX} \\
& U \\
U \ar[ru]_-{\one_U} \ar[rr]_-{\Delta_U} & & 
U\times U \ar[lu]^-{j^*_UU}
\enddiagram}\qquad
$$
where $\Gamma_{\!\phi}$ is the graph of $\phi$, and $s^*\cD$ is
isomorphic to the cartesian diagram (in $T'$) :
$$
\cE
\quad : \quad
{\spreaddiagramcolumns{+30pt}\diagram
s^*\phi \ar[r]^-{\omega^\dagger_X\circ s^*\Gamma_\phi}
\ar[d]_{u_{s^*\one_U}\circ s^*(\phi/U)} & f^*X \ar[d]^{f^*\phi} \\
1_{T'} \ar[r]^-{\sigma(s,\omega)} & f^*U.
\enddiagram}
\qquad\qquad\qquad
$$
This shows that $s^*$ -- and therefore also $s$ -- is determined,
up to isomorphism, by $\sigma(s,\omega)$.

(ii)\ \
Conversely, if $\sigma\in\Gamma(T',f^*U)$ is any global section,
we may define a functor $s^*:T\!/U\to T'$ by means of the cartesian
diagram $\cE$ : this amounts to choosing a representative in $T'$
for the fibre product $1_{T'}\times_{(\sigma,f^*\phi)}f^*X$, for every
object $\phi:X\to U$ of $T/U$. We claim that there exists an
isomorphism of functors
$$
s^*\circ j_U^*\isom f^*.
$$
Indeed, recall that $j^*_UX=(q_X:X\times U\to U)$, where $X\times U$
is a representative of the product of $X$ and $U$ in $T$, and
$q_X$ is the natural projection; let also $p_X:X\times U\to X$
be the other projection. Likewise, pick a representative
$f^*X\times f^*U$ for the fibre product of $f^*X$ and $f^*U$ in
$T'$, and let $p_{f^*X}:f^*X\times f^*U\to f^*X$ and
$q_{f^*X}:f^*X\times f^*U\to f^*U$ be the projections.
Since $f^*$ is left exact, there exists a unique isomorphism
$\lambda_X:f^*X\times f^*U\isom f^*(X\times U)$ in $T'$ such that
$f^*(p_X)\circ\lambda_X=p_{f^*X}$ and $f^*(q_X)\circ\lambda_X=q_{f^*X}$.
We deduce a commutative diagram :
$$
\xymatrix@C+50pt{
f^*X \ar[r]^-{\one_{f^*X}\times(\sigma\circ u_{f^*X})} \ar[d]_{u_{f^*X}}
& f^*X\times f^*U \ar[d]^{q_{f^*X}} \ar[r]^-\lambda &
f^*(X\times U) \ar[ld]^{f^*j^*_UX} \\
1_{T'} \ar[r]^-\sigma & f^*U
}$$
whose square subdiagram is cartesian (the details are left
to the reader). The assertion follows straightforwardly.
Let us also define a functor $t:T'\to T'\!/f^*U$, by the rule :
$$
tY:=\sigma\circ u_Y
\qquad
\text{for every $Y\in\Ob(T')$}.
$$
By inspecting the diagram $\cE$, we deduce natural bijections :
$$
\Hom_{T'}(Y,s^*\phi)\isom\Hom_{T'\!/f^*U}(tY,f^*\phi)
\qquad
\text{for every $Y\in\Ob(T')$ and $\phi\in\Ob(T\!/U)$.}
$$
Since $f^*$ is exact, it follows easily that $s^*$ is left
exact (indeed, $s^*$ commutes with all the limits with
which $f^*$ commutes). Moreover, since all colimits are universal
in $T$ (see \eqref{rem_rep-and-sheafify}(i)), it is easily seen
that $s^*$ commutes with all colimits. Then, by proposition
\ref{prop_half-dominates}(i.b), the functor $s^*$ determines a
morphism of topoi $s:T'\to T\!/U$, unique up to unique
isomorphism, with an isomorphism $f\isom j_U\circ s$.

(iii)\ \
Next, let $s,s':T'\to T/U$ be two morphisms of topoi,
$\omega:f\isom j_U\circ s$ and $\omega':f\isom j_U\circ s'$
two isomorphisms, and $\beta:s\Rightarrow s'$ any natural
transformation such that
\set\begin{equation}\label{eq_ischia-ok}
(j_U*\beta)\odot\omega=\omega'.
\end{equation}
By remark \ref{rem_adjoint-transf}(iv,vi), it follows that
$\omega'^\dagger=\omega^\dagger\odot(\beta^\dagger*j_U^*)$.
We may then compute :
$$
\begin{aligned}
\sigma(s',\omega')=
\omega'^\dagger_U\circ s'^*(\Delta_U)\circ u^{-1}_{s'^*\one_U}
=&\,\omega^\dagger_U\circ\beta^\dagger_{j_U^*U}\circ s'^*(\Delta_U)
\circ u^{-1}_{s'^*\one_U} \\
=&\,\omega^\dagger_U\circ s^*(\Delta_U)\circ
\beta^\dagger_{\one_U}\circ u^{-1}_{s'^*\one_U} \\
=&\,\sigma(s,\omega).
\end{aligned}
$$
By considering the cartesian diagram $\cE$ of (i), we then
deduce that for every object $\phi:X\to U$ of $T/U$ there
exists a unique isomorphism
$$
\lambda_\phi:s'^*\phi\isom s^*\phi
\qquad\text{such that}\qquad
\omega^\dagger_X\circ s^*\Gamma_\phi\circ\lambda_\phi=
\omega'^\dagger_X\circ s'^*\Gamma_\phi.
$$
Since $\omega^\dagger_X$ is an isomorphism, the latter
identity yields :
$$
s^*\Gamma_\phi\circ\lambda_\phi=s^*\Gamma_\phi\circ\beta^\dagger_\phi.
$$
However, $\Gamma_\phi$ is a monomorphism, hence the same holds
for $s^*\Gamma_\phi$ (proposition \ref{prop_was-also-cofinal}(i)),
whence $\lambda_\phi=\beta^\dagger_\phi$. Hence, $\beta$ is
{\em the unique isomorphism of functors} $s\isom s'$ verifying
\eqref{eq_ischia-ok}.

(iv)\ \
Lastly, consider a morphism $\psi:U\to V$ of $T$, and
$\sigma\in\Gamma(T',f^*U)$. Let $s^U:T/U\to T'$ be the morphism
of topoi deduced from $\sigma$ as in (ii), with its associated
isomorphism of functors $\omega^U:s^{U*}\circ j_U^*\isom f^*$.
Likewise, let $s^V:T/V\to T'$ be the morphism of topoi and
$\omega^V:s^{V*}\circ j_V^*\isom f^*$ the isomorphism deduced,
again as in (ii), from $f^*(\psi)\circ\sigma\in\Gamma(T',f^*V)$.
A direct inspection of the constructions yields a natural
isomorphism of functors
$$
\beta:s^{V*}\isom s^{U*}\circ j_\psi^*
\qquad\text{such that}\qquad
\omega^V=
\omega^U\odot(s^{U*}*\gamma^{\dagger\,-1}_{\psi,V})\odot(\beta*j^*_V)
$$
where $\gamma_{\psi,V}:j_V\circ j_\psi\isom j_U$ is the
coherence constraint of the fibred topos $T/-$ of example
\ref{ex_localization-topos}(ii). In view of (iii), it follows
easily that for every pair $(s,\omega)$ as in (i) we have :
$$
\sigma(j_\psi\circ s,(\gamma^{-1}_{\psi,V}*s)\odot\omega)=
f^*(\psi)\circ\sigma(s,\omega)
$$
\end{remark}

\begin{definition}\label{def_T-point}
Let $T$ be a topos.

(i)\ \
A {\em point\/} of $T$ (or a {\em $T$-point}) is a morphism of topoi
$\Set\to T$. If $\xi=(\xi^*,\xi_*)$ is a point, and $F$ is any
object of $T$, the set $\xi^*F$ is usually denoted by $F_\xi$.

(ii)\ \
If $\xi$ is a $T$-point, and $f:T\to S$ is any morphism of topoi,
we denote by $f(\xi)$ the $S$-point $f\circ\xi$. If $\beta:\xi\to\xi'$
is any morphism of $T$-points, we let likewise $f(\beta):=f_**\beta$.

(iii)\ \
A {\em neighborhood\/} of $\xi$ is a pair $(U,a)$, where
$U\in\Ob(T)$, and $a\in U_\xi$. A morphism of neighborhoods
$(U,a)\to(U',a')$ is a morphism $f:U\to U'$ in $T$ such that
$f_\xi(a)=a'$. The category of all neighborhoods of $\xi$
shall be denoted $\bNbd(\xi)$.

(iv)\ \
A family $(f_\lambda:T_\lambda\to T~|~\lambda\in\Lambda)$ of
morphisms of topoi is {\em conservative}, if the functor
$$
T\to\prod_{\lambda\in\Lambda}T_\lambda
\qquad
F\mapsto(f^*_\lambda F~|~\lambda\in\Lambda)
$$
is conservative (definition \ref{def_equivalence}(ii)) (the product
$\prod_{\lambda\in\Lambda}T_\lambda$ is formed in a sufficiently large
universe containing $\sU$). We say that $T$ {\em has enough
points}, if $T$ has a conservative set of points.

(v)\ \
Consider a graph $\Gamma$ (see definition \ref{def_graph}(i)) and
a system $(C_S~|~S\in\Ob(\Topos))$, where $C_S$ is a given set of
morphisms of graphs $\Gamma\to S$, for every topos $S$, such that
$f^*\phi\in C_{S'}$, for every $\phi\in C_S$ and every morphism
$f:S'\to S$ of topoi.
Let $T_\bullet:=(T_\lambda\xrightarrow{f_\lambda}T~|~\lambda\in\Lambda)$
be a conservative family of morphisms of topoi, and $\bP(S,\phi)$
a property defined on elements $\phi\in C_S$, for every
topos $S$; we say that $\bP$ {\em can be checked on} $T_\bullet$,
if for every $\phi\in C_T$ we have :
\begin{itemize}
\item
$\bP(T,\phi)$ holds if and only if
$\bP(T_\lambda,f^*_\lambda\circ\phi)$ holds for every $\lambda\in\Lambda$.
\end{itemize}
Moreover, we say that $\bP$ {\em can be checked on stalks}, if for
every topos $S$ and every conservative set $\Omega$ of $S$-points,
the property $\bP$ can be checked on $\Omega$.
\end{definition}

\begin{example}\label{ex_covers-are-conserv}
Let $T$ be a topos, and $(U_\lambda\to 1_T~|~\lambda\in\Lambda)$
any family of morphisms of $T$ covering the final object $1_T$
(for the canonical topology of $T$); then the induced family
of morphisms of topoi
$(j_{U_\lambda}:T/U_\lambda\to T~|~\lambda\in\Lambda)$ is conservative.
For the proof, we easily reduce to the case where $\Lambda$ is
a small set, by lemma \ref{lem_U-site}; then it suffices to
recall that, by remark \ref{rem_summarized}(i), every $X\in\Ob(T)$
is the coequalizer of the induced pair of morphisms of $T$
$$
\xymatrix{\coprod_{\lambda,\mu\in\Lambda}U_\lambda\times U_\mu\times X
\ar@<+.5ex>[r] \ar@<-.5ex>[r] & \coprod_{\lambda\in\Lambda}U_\lambda\times X}
$$
(the details are left to the reader).
\end{example}

\begin{remark}\label{rem_Nbd-as-Fib}
(i)\ \
With the notation of definition \ref{def_T-point},
notice that
$$
\bNbd(\xi)^o=\cFib(\xi^*)
$$
where the fibration $\pi:\cFib(\xi^*)\to T^o$ is defined
as in \eqref{subsec_fibred-cats-II}.

(ii)\ \
Moreover, the category $\bNbd(\xi)$ is finitely complete.
Indeed, let $\Lambda$ be any finite category, and
$F:\Lambda\to\bNbd(\xi)$ any functor; say that
$F\lambda=(U_\lambda,a_\lambda)$ for every $\lambda\in\Ob(\Lambda)$.
Then it is easily seen that the limit of $F$ is represented
by $(U,a)$, where $U\in\Ob(T)$ represents the limit of the
functor $\pi^o\circ F:\Lambda\to T$, and
$a\in U_\xi=\lim_\Lambda\xi^*\circ\pi^o\circ F$ is the unique
element whose image under the induced projection
$U_\xi\to U_{\lambda,\xi}$ agrees with $a_\lambda$, for every
$\lambda\in\Ob(\Lambda)$.

(iii)\ \
Let $\xi$ be a point of the topos $T$; from (ii) it follows
that the category $\bNbd(\xi)$ is cofiltered.

(iv)\ \
In view of (i) and \eqref{subsec_category-of-elements}, for
every $F\in\Ob(T)$ there is a natural cocone :
$$
\xymatrix@C+30pt{
\bNbd(\xi)^o \rtwocell^{\Hom_T(\pi,F)}_{c_{F_\xi}}{\ \ \ \tau_{\xi,F}} &
\Set
}$$
assigning to every neighborhood $(U,a)$ of $\xi$,
the map $\tau_{\xi,F}(U,a):\Hom_T(U,F)\to F_\xi$ such that
$s\mapsto s_\xi(a)$ for every $s:U\to F$. Then lemma
\ref{lem_lable} says that $\tau_{\xi,F}$ induces a natural
bijection
\set\begin{equation}\label{eq_deliver-evil}
\colim_{\bNbd(\xi)^o}\Hom_T(\pi,F)\isom F_\xi
\qquad
\text{for every $F\in\Ob(T)$}.
\end{equation}

(v)\ \
It follows from remark \ref{rem-sections-of-f}(iii), that a
neighborhood $(U,a)$ of $\xi$ is the same as the datum of an
isomorphism class of a point $\xi_U$ of the topos $T\!/U$ which
lifts $\xi$, {\em i.e.} such that $\xi\simeq j_U(\xi_U)$.
Moreover, say that $\xi_{U'}$ is another lifting of $\xi$,
corresponding to a neighborhood $(U',a')$ of $\xi$; then, by
inspecting the constructions of remark
\ref{rem-sections-of-f}(i,ii) we see that, under this
identification, a morphism $(U,a)\to(U',a')$ of neighborhoods
of $\xi$ corresponds to the datum of :
$$
\text{a morphism $\phi:U\to U'$ in $T$ and an isomorphism
$j_\phi(\xi_U)\isom \xi_{U'}$ of $T$-points}
$$
where $j_\phi$ is the morphism of topoi $T/U\to T/U'$ induced
by $\phi$.
\end{remark}

\begin{proposition}\label{prop_formal-chenck-stalks}
Let $T$ be a topos with enough points; we have :

{\em (i)}\ \
For every object $U$ of\/ $T$, the topos $T/U$ has enough points.

{\em (ii)}\ \
More precisely, let $\Omega$ be a conservative set of\/ $T$-points,
and denote by $j^{-1}_U\Omega$ the set of all $T/U$-points $\xi_U$
such that there exists $\xi\in\Omega$ with an isomorphism
$\xi\isom j_U(\xi_U)$; then $j^{-1}_U\Omega$ is a conservative
set of\/ $T/U$-points.
\end{proposition}
\begin{proof} Let $\phi:X\to Y$ be a morphism in $T/U$ such that
$\phi_\xi$ is an epimorphism for every $\xi\in j^{-1}_U\Omega$; by
proposition \ref{prop_conser-nd-left-exact}(ii), it suffices to
show that $\phi$ is an epimorphism, and the latter will follow,
if we show that the same holds for $j_{U!}\phi$.
Thus suppose by way of contradiction, that $j_{U!}\phi$
is not an epimorphism; then there exist $\xi\in\Omega$ and
$y\in(j_{U!}Y)_\xi$, such that $y$ does not lie in the image of
$(j_{U!}\phi)_\xi$. Let $a:=(j_{U!}\pi_Y)_\xi(y)\in U_\xi$, where
$\pi_Y:Y\to 1_{T\!/U}=\one_U$ is the unique morphism in
$T/U$. By remark \ref{rem_Nbd-as-Fib}(v), we may find a
lifting $(\xi_U:\Set\to T/U,\omega_U)$ of $\xi$ such that
\set\begin{equation}\label{eq_condition}
\sigma(\xi_U,\omega_U)=a
\end{equation}
where $\sigma(\xi_U,\omega_U)$ is defined as in remark
\ref{rem-sections-of-f}(i).
After replacing $\xi$ by $j_U(\xi_U)$, we may assume that
$\omega_U=\one_\xi$, in which case \eqref{eq_condition} means
that $a=\xi_U^*\Delta_U$, where
$\Delta_U:1_{T\!/U}\to j^*_Uj_{U!}1_{T\!/U}$ is the unit
of adjunction ({\em i.e} the diagonal $U\to U\times U$).
We have a commutative diagram of sets :
$$
\xymatrix{
\xi_U^*X \ar[rr]^-{\xi^*_U\phi} \ar[d]_{\xi_U^*\Delta_X}
& & \xi_U^*Y \ar[d]^{\xi_U^*\Delta_Y} \\
(j_{U!}X)_\xi \ar[rr]^-{(j_{U!}\phi)_\xi} & &
(j_{U!}Y)_\xi=\xi_U^*(j_U^*j_{U!}Y)
}$$
where $\Delta_X:X\to j^*_Uj_{U!}X$ is the unit of adjunction,
and likewise for $\Delta_Y$. More plainly,
$\Delta_Y:Y\times 1_{T\!/U}\to Y\times(j^*j_{U!}1_{T\!/U})$
is the product $\one_Y\times\Delta_U$, and under this
identification, $\xi_U^*\Delta_Y$ is the mapping
$\xi^*_UY\to\xi^*_UY\times(j^*j_{U!}1_{T\!/U})$ given by
the rule : $z\mapsto(z,a)$ for every $z\in\xi_U^*Y$.
Especially, we see that $y$ lies in the image of
$\xi^*_U\Delta_Y$. But by assumption, the map $\xi^*_U\phi$
is surjective, hence $y$ lies in the image of $(j_{U!}\xi)_\xi$,
a contradiction.
\end{proof}

\sset\subsubsection{}\label{subsec_blift}
Remark \ref{rem_Nbd-as-Fib}(v) prompts the following construction.
Let us consider the fibred topos
$$
T/-:T\to\Topos
\qquad
U\mapsto T/U
\qquad
(U'\xrightarrow{\phi}U)\mapsto(T/U'\xrightarrow{j_\phi}T/U) 
$$
of example \ref{ex_localization-topos}(ii). Denote also by $1_T$
a fixed final object of $T$, and for every $U\in\Ob(T)$ let
$\phi_U:U\to 1_T$ be the unique morphism in $T$; the source
functor yields a natural isomorphism $T/1_T\isom T$ that
identifies $j_{\phi_U}:T/U\to T/1_T$ with $j_U:T/U\to T$.
Moreover, the rule : $U\mapsto\phi_U$ clearly defines a
natural transformation
$$
\phi_\bullet:\one_T\Rightarrow c_{1_T}
$$
where $c_{1_T}:T\to T$ is the constant functor with value $1_T$.
For a sufficiently large universe $\sV$, example
\ref{ex_first-from-Cats-to-PsFuns}(i) yields a strict
pseudo-functor $H_\Topos:\Topos^o\times\Topos\to\sV\tdu\bCat$
which restricts to a strict pseudo-functor
$$
\Pt:\Topos\to\sV\tdu\bCat
\qquad
S\mapsto\Pt(S):=\Hom_\Topos(\Set,S)
$$
assigning to every topos $S$ the category of $S$-points.
We consider the composition
$$
\Pt_T:=\Pt\circ(T/-):T\to\sV\tdu\bCat
\qquad
U\mapsto\Pt(T/U)
$$
as well as the pseudo-natural transformation
$$
\Pt_T*\phi_\bullet:\Pt_T\Rightarrow\sF_{\Pt(T)}
$$
from $\Pt_T$ to the constant pseudo-functor $\sF_{\Pt(T)}$
with value $\Pt(T)$. There follows a cartesian functor of
fibrations over $T^o$ :
$$
\cFib(\Pt_T*\phi_\bullet):\cFib(\Pt_T)\to T^o\times\Pt(T).
$$
Let also $\sigma_\xi:T^o\to T^o\times\Pt(T)$ be the cartesian
section such that $\sigma_\xi(U):=(U,\xi)$ for every $U\in\Ob(T)$.
Lastly, we consider the $2$-cartesian diagram in the category
$\Fib(T^o)$ :
$$
\xymatrix{
\cL(\xi):=\displaystyle{\cFib(\Pt_T)\mathop{\times}^2_{T^o\times\Pt(T)}T^o}
\ar[r] \ar[d] & \cFib(\Pt_T) \ar[d]^{\cFib(\Pt_T*\phi_\bullet)} \\
T^o \ar[r]^{\sigma_\xi} & T^o\times\Pt(T) 
}$$
and we set
$$
\bLift(\xi):=\cL(\xi)^o.
$$
Taking into account theorem \ref{th_fib-limits-are-fibrewise}(ii),
we see that for every $U\in\Ob(T)$, the fibre $\cL(\xi)_U$ of the
fibration $\cL(\xi)\to T^o$ is the category whose objects are
the pairs $(\xi_U,\omega_U)$, where $\xi_U$ is a $T/U$-point,
and $\omega_U:\xi\isom j_U(\xi_U)$ is an isomorphism of
$T$-points. The morphisms $(\xi_U,\omega_U)\to(\xi'_U,\omega'_U)$
are the morphisms of $T/U$-points
$$
\beta:\xi_U\to\xi'_U
\qquad\text{such that}\qquad
\omega_{U'}=j_U(\beta)\odot\omega_U
$$
(with the obvious composition law). Hence, the objects
of $\bLift(\xi)$ are the triples $(U,\xi_U,\omega_U)$
with $U\in\Ob(T)$ and $(\xi_U,\omega_U)\in\Ob(\cL(\xi)_U)$.
The morphisms $(U,\xi_U,\omega_U)\to(V,\xi_V,\omega_V)$
are the pairs $(\phi,\beta)$, where $\phi:U\to V$ is a
morphism in $T$, and $\beta:\xi_V\to j_\phi(\xi_U)$
is a morphism of $T\!/V$-points, such that
$$
(\gamma_{\phi,V}*\xi_U)\odot j_V(\beta)\odot\omega_V=\omega_U
$$
where $\gamma_{\phi,V}:j_V\circ j_\phi\isom j_U$ is the coherence
constraint of $T/-$. Then, remark \ref{rem-sections-of-f} shows
that $\cL(\xi)$ is a fibration in groupoids, and we have an
equivalence of categories :
\set\begin{equation}\label{eq_same-neighbor}
N_\xi:\bLift(\xi)\isom\bNbd(\xi)
\qquad
(U,\xi_U,\omega_U)\mapsto(U,\sigma(\xi_U,\omega_U))
\end{equation}
where $\sigma(\xi_U,\omega_U)\in U_\xi$ is defined as in
remark \ref{rem-sections-of-f}(i).

\sset\subsubsection{}\label{subsec_to-every-set}
Notice next that, since $\cL(\xi)$ is a fibration in groupoids,
the natural projection $\cL(\xi)\to\cFib(\Pt_T)$ factors
through the inclusion pseudo-functor
$\cFib(\Pt_T^\times)=\cFib(\Pt_T)^\times\to\cFib(\Pt_T)$ (notation
of remark \ref{rem_groupoids}(iii)). To every small set $X$, we
wish now to attach a functor
$$
\Xi^X:\cFib(\Pt^\times)\to\Set
\qquad
(\zeta:\Set\to S)\mapsto\zeta^*\zeta_*X.
$$
To this aim, according to remark \ref{rem_functors-from-lax}(i),
it will suffice to exhibit a lax-natural transformation
$$
\alpha^X:{}^o\Pt^\times\Rightarrow{}^o\sF_\Set
$$
where $\sF_\Set$ is the constant pseudo-functor
$\Topos\to\sV\tdu\bCat$ with value $\Set$. We need
therefore to give for every $S\in\Ob(\Topos)$ a functor
$\alpha^X_S:\Pt(S)^\times\to\Set$, together with a coherence
constraint for the resulting system
$(\alpha^X_S~|~S\in\Ob(\Topos))$. Now, let us set
$$
\alpha^X_S(\xi):=\zeta^*\zeta_*X
\qquad
\text{for every $\zeta\in\Ob(\Pt(S))$}.
$$
Next, let $\beta:\zeta\isom\zeta'$ be any morphism in
$\Pt(S)^\times$; {\em i.e.} $\beta:\zeta_*\isom\zeta'_*$ is an
isomorphism of functors. Then we define
$\alpha^X_S(\beta):\zeta^*\zeta_*X\to\zeta'^*\zeta'_*X$ as
the composition
$$
\zeta^*\zeta_*X\xrightarrow{(\zeta^**\beta)_X}\zeta^*\zeta'_*X
\xrightarrow{(\beta^\dagger*\zeta'_*)^{-1}_X}\zeta'^*\zeta'_*X
$$
where $\beta^\dagger:\zeta'^*\isom\zeta^*$ is the transpose
of the natural transformation $\beta$ (see remark
\ref{rem_adjoint-transf}(ii)). From remark
\ref{rem_adjoint-transf}(vi) we see already that
$\alpha^X_S(\one_\zeta)=\one_{\zeta^*\zeta_*X}$. Next, let
$\beta':\zeta'\to\zeta''$ be another morphism of $\Pt(S)^\times$;
we compute :
$$
\begin{aligned}
\alpha^X_S(\beta')\circ\alpha^X_S(\beta)&\,=
((\beta'^\dagger*\zeta''_*)^{-1}\odot(\zeta'^**\beta')\odot
(\beta^\dagger*\zeta'_*)^{-1}\odot(\zeta^**\beta))_X \\
&\,=((\beta'^\dagger*\zeta''_*)^{-1}\odot
(\beta^\dagger*\zeta''_*)^{-1}\odot(\zeta^**\beta')
\odot(\zeta^**\beta))_X \\
&\,=(((\beta'\odot\beta)^\dagger*\zeta''_*)^{-1}\odot
(\zeta^**(\beta'\odot\beta)))_X \\
&\,=\alpha^X_S(\beta'\odot\beta)
\end{aligned}
$$
where the second equality follows from
\eqref{eq_Godement-functor}, and the third follows
from remark \ref{rem_adjoint-transf}(vi).
Next, let $f:S\to S'$ be any morphism of topoi, and
$\eps^f:f^*\circ f_*\Rightarrow\one_S$ the counit of the
adjoint pair $(f^*,f_*)$; our candidate coherence constraint
for $\alpha^X_\bullet$ assigns to every $S$-point $\zeta$ the map
$$
\tau^X_{f,\zeta}:f(\zeta)^*f(\zeta)_*X=
\zeta^*\circ f^*\circ f_*\circ\zeta_*X
\xrightarrow{(\zeta^**\eps^f*\zeta_*)_X}\zeta^*\zeta_*X
$$
(see example \ref{ex_from-Cats-to-PsFuns}(ii)). We need
to check that the rule $\zeta\mapsto\tau^X_{f,\zeta}$ yields
a natural transformation
$$
\tau^X_f:\alpha^X_{S'}\circ\Pt(f)^\times\Rightarrow\alpha^X_S.
$$
Thus, let $\beta:\zeta_*\isom\zeta'_*$ be any morphism of
$\Pt(S)^\times$; the assertion comes down to showing the
commutativity of the diagram :
$$
\xymatrix@C+60pt{ \zeta^*f^*f_*\zeta_*X
\ar[r]^-{(\zeta^*f^*f_**\beta)_X} \ar[d]_{(\zeta^**\eps^f*\zeta_*)_X} &
\zeta^*f^*f_*\zeta'_*X \ar[r]^-{((f_**\beta)^\dagger*f_*\zeta'_*)^{-1}_X}
\ar[d]_{(\zeta^**\eps^f*\zeta'_*)_X} &
\zeta'^*f^*f_*\zeta'_*X \ar[d]^{(\zeta'^**\eps^f*\zeta'_*)_X} \\
\zeta^*\zeta_*X \ar[r]^-{(\zeta^**\beta)_X} &
\zeta^*\zeta'_*X \ar[r]^-{(\beta^\dagger*\zeta'_*)^{-1}_X} &
\zeta'^*\zeta'_*X
}$$
However, the commutativity of the left square subdiagram
follows from the naturality of $\zeta^**\eps^f$. Next,
recall that $(f_**\beta)^\dagger=\beta^\dagger*f^*$ (remark
\ref{rem_adjoint-transf}(vi)); then the commutativity
of the right square subdiagram follows by applying 
\eqref{eq_Godement-functor} to the natural transformations
$\beta^\dagger$ and $\eps^f*\zeta'_*$.

Lastly, we have to check the coherence axioms for $\tau^X_\bullet$.
Since ${}^o\Pt^\times$ and ${}^o\sF_\Set$ are both strict
pseudo-functors, we come down to showing that
$$
\tau^X_{\one_S}=\one_{\alpha^X_S}
\qquad\text{and}\qquad
\tau^X_g\boxminus\tau^X_f=\tau^X_{g\circ f}
$$
for every topos $S$ and every pair of morphisms of
topoi $S\xrightarrow{f}S'\xrightarrow{g}S''$ (see
remark \ref{rem_pseudo-natural}(i)). These follow
easily by direct inspection. This completes the
construction of $\alpha^X$, and thus also of
$\Xi^X:=\cFib(\alpha^X)$. After composing with the
projection $\cL(\xi)\to\cFib(\Pt^\times_T)\to\cFib(\Pt^\times)$,
we deduce a functor
$$
\Xi^{X,\xi}:\bLift(\xi)^o\to\Set
\qquad
(U,\xi_U,\omega_U)\mapsto\xi_U^*\xi_{U*}X.
$$

\begin{lemma}\label{lem_super-cazzuto}
{\em(i)}\ \
For every $(S,\zeta)\in\Ob(\cFib(\Pt^\times))$ and every set
$X$, let $\eps^{(S,\zeta)}_X:\zeta^*\zeta_*X\to X$ be the
counit of adjunction. Then the rule :
$(S,\zeta)\mapsto\eps^{(S,\zeta)}_X$ defines a natural cocone
$$
\Xi^X\Rightarrow c_X.
$$

{\em(ii)}\ \
The cocone of\/ {\em (i)} induces a natural bijection :
$$
\colim_{\bLift(\xi)^o}\Xi^{X,\xi}\isom X
\qquad
\text{for every (small) set $X$.}
$$
\end{lemma}
\begin{proof}(i): Let $\nu^X:{}^o\Pt^\times\Rightarrow{}^o\sF_\Set$
be the strict pseudo-natural transformation such that
$$
\nu^X_S:\Pt(S)^\times\to\Set
$$
is the constant functor with value $X$, for every topos $S$
({\em i.e.} $\nu^X_S(\zeta):=X$ for every $S$-point $\zeta$,
and $\nu^X_S(\beta)=\one_X$ for every invertible morphism
$\beta$ of $S$-points). For every such $S$, let also
$$
\eps^{(S,\bullet)}_X:\Ob(\Pt(S))\to\rMorph(\Set)
$$
be the map that associates with every $S$-point $\zeta$,
the map of sets $\eps^{(S,\zeta)}_X$. We notice :

\begin{claim} The rule : $S\mapsto\eps^{(S,\bullet)}_X$ for
every topos $S$, defines a modification :
$$
\eps^{(\bullet,\bullet)}_X:\nu^X\leadsto\alpha^X.
$$
\end{claim}
\begin{pfclaim} First, we need to check that $\eps^{(S,\bullet)}_X$
is a natural transformation $\alpha^X_S\Rightarrow\nu^X_S$, for
every topos $S$. The latter assertion amounts to the
commutativity of the diagram :
$$
\xymatrix@C+20pt{ \zeta^*\zeta_*X \ar[r]^-{\eps^{(S,\zeta)}_X}
\ar[d]_{\alpha^X_S(\beta)} & X \ddouble \\
\zeta'^*\zeta'_*X \ar[r]^-{\eps^{(S,\zeta')}_X} & X
}$$
for every isomorphism $\beta:\zeta\isom\zeta'$ of $S$-points.
Since $(\beta^\dagger*\zeta'_*)^{-1}_X=((\beta^{-1})^\dagger*\zeta'_*)_X$
(remark \ref{rem_adjoint-transf}(iv)), this is in turn the same
as the identity :
$$
\eps^{(S,\zeta')}_X\circ((\beta^{-1})^\dagger*\zeta'_*)_X=
\eps^{(S,\zeta)}_X\circ(\zeta^**\beta^{-1})_X
$$
which follows from remark \ref{rem_adjoint-transf}(iii).
Next, we need to check the compatibility condition for
$\eps^{(S,\bullet)}_X$, which comes down to the identity :
$$
\eps^{(S,\zeta)}_X\odot\tau^X_f=\eps^{(S,f(\zeta))}_X
$$
for every morphism of topoi $f:S\to S'$ and every $S$-point
$\zeta$, where $\tau^X_f$ denotes the coherence constraint
of $\alpha^X$. The latter follows by direct inspection.
\end{pfclaim}

By remark \ref{rem_functors-from-lax}(ii), the sought cocone
is then
$\cFib(\eps^{(\bullet,\bullet)}_X):\Xi^X\Rightarrow\cFib(\nu^X)=c_X$.

(ii): Let $\pi:\cFib(\Pt^\times)\to\Topos^o$ be the fibration
arising from the strict pseudo-functor $\Pt^\times$. We let $\cC$
be the fibre product in the cartesian square of categories :
$$
\xymatrix{ \sMorph(\cFib(\Pt^\times)) \ar[r]^-p \ar[rd]_{\sMorph(\pi)} &
\cC \ar[r]^-{\ss'} \ar[d]_{\pi'} & \cFib(\Pt^\times) \ar[d]^\pi \\
& \sMorph(\Topos^o) \ar[r]^-\ss & \Topos^o
}$$
where $\ss$ and $\ss'\circ p$ are the source functors. By remark
\ref{rem_distinguished-cleavage}(ii), we have a natural
identification
$$
\cC\isom\cFib(\Pt^\times\circ\ss^o)
$$
and recall that $\sMorph(\pi)$ is a fibration as well,
by example \ref{ex_Morph-is-fibred}. To every morphism of
topoi $f:S'\to S$, we attach the functor
$$
\Sigma^X_f:\pi'^{-1}(f)=\Pt(S)^\times\to\Set
\qquad
\zeta\mapsto\Gamma(S',f^*\zeta_*X)
$$
where $\Gamma(S',-):S'\to\Set$ is defined as in
\eqref{subsec_glob-sections}, after fixing a final object
$1_{S'}$ of $S'$. To every isomorphism
$\alpha:\zeta_1\isom\zeta_2$ of $S$-points, the functor
$\Sigma^X_f$ assigns the map
$$
\Sigma^X_f(\alpha):=\Gamma(S',f^*\alpha_X):
\Gamma(S',f^*\zeta_{1*}X)\to\Gamma(S',f^*\zeta_{2*}X).
$$
Next, notice that if $(S'_1\xrightarrow{f_1}S_1)^o$ and
$(S'_2\xrightarrow{f_2}S_2)^o$ are two objects of $\sMorph(\Topos^o)$,
then a morphism $f^o_1\to f^o_2$ in this category is a commutative
diagram of morphisms of topoi :
$$
\cD_{g,g'} \qquad :\qquad
{\diagram S'_2 \ar[r]^-{f_2} \ar[d]_{g'} & S_2 \ar[d]^g \\
S'_1 \ar[r]^-{f_1} & S_1.
\enddiagram}\qquad\qquad\qquad
$$
The (split) cleavage $\Pt^\times\circ\ss^o$ of $\cC$
associates with $\cD_{g,g'}$ the functor
$$
\Pt(g)^\times:\Pt(S_2)^\times\to\Pt(S_1)^\times
\qquad
\zeta\mapsto g(\zeta)
$$
that assigns to every morphism $\alpha$ of $\Pt(S_2)^\times$
the morphism $g(\alpha)$ of $\Pt(S_1)^\times$. To the diagram
$\cD_{g,g'}$ and every object $\zeta$ of $\Pt(S_2)^\times$ we
attach as well the map :
$$
\tau^{g,g'}_\zeta:
\Gamma(S'_1,f_1^*g_*\zeta_*X)
\xrightarrow{\ \Gamma(S'_1,\Upsilon(\cD_{g,g'})*\zeta_*)_X\ }
\Gamma(S'_1,g'_*f_2^*\zeta_*X)
\xrightarrow{\omega^{g'}_{f_2^*\zeta_*X}}\Gamma(S'_2,f_2^*\zeta_*X)
$$
where $\Upsilon(\cD_{g,g'})$ is the base change transformation
associated with the diagram $\cD_{g,g'}$, viewed as an oriented
square diagram of links, whose orientation is the identity
(see \eqref{subsec_base-change-map}). The second map in
this composition is induced by the unique isomorphism of
morphisms of topoi $\omega^{g'}:\Gamma_{S'_1}\circ g'_*\isom\Gamma_{S'_2}$
(proposition \ref{prop_glob-sections}(ii)).

\begin{claim}\label{cl_patience}
(i)\ \
For every diagram $\cD_{g,g'}$, the rule :
$\zeta\mapsto\tau^{g,g'}_\zeta$ defines a natural transformation
$$
\xymatrix@C+20pt{ \Pt(S_2)^\times \ar[r]^-{\Sigma^X_{f_2}}
\ar[d]_{\Pt(g)^\times} \drtwocell\omit{^\tau^{g,g'}\ \ \ \ } &
\Set \ddouble \\
\Pt(S_1)^\times \ar[r]_-{\Sigma^X_{f_1}} & \Set.
}$$

(ii)\ \
The rule : $\cD_{g,g'}\mapsto\tau^{g,g'}$ provides the coherence
constraint for a lax-natural transformation
$$
\Sigma^X_\bullet:{}^o(\Pt^\times\circ\ss^o)\Rightarrow{}^o\sF_\Set.
$$
\end{claim}
\begin{pfclaim}(i): For every object $(S'\xrightarrow{f}S)^o$
of $\sMorph(\Topos^o)$, consider the functor
$$
\sigma_f:\Pt(S)^\times\to S
\qquad
\zeta\mapsto\zeta_*X
$$
that assigns to every isomorphism $\alpha$ as in the
foregoing, the morphism $\alpha_X$ of $S$. With this notation,
$\tau^{g,g'}$ is the composition of the following oriented
squares :
$$
\xymatrix@C+60pt{ \Pt(S_2)^\times \ar[r]^-{\sigma_{f_2}}
\ar[d]_{\Pt(g)^\times}
\drtwocell\omit{^\one_{g_*\circ\sigma_{f_2}}\ \ \ \ \ \ \ \ }
& S_2 \ar[r]^-{f^*_2} \ar[d]_{g_*}
\drtwocell\omit{^\Upsilon(\cD_{g,g'})\ \ \ \ \ \ \ \ \ } &
S'_2 \ar[r]^-{\Gamma_{S'_2}} \ar[d]^{g'_*}
\drtwocell\omit{^\omega^{g'}\ \ \ } & \Set \ddouble \\
\Pt(S_1)^\times \ar[r]_-{\sigma_{f_1}} &
S_1 \ar[r]_-{f^*_1} & S'_2 \ar[r]_-{\Gamma_{S'_1}} & \Set.
}$$

(ii): Clearly $\tau^{\one_S,\one_{S'}}=\one_{\Sigma^X_f}$, for every
morphism of topoi $f:S'\to S$. According to remark
\ref{rem_pseudo-natural}(i), it remains to check the identity
$$
\tau^{g,g'}\boxminus\tau^{h,h'}=\tau^{g\circ h,g'\circ h'}
$$
for every composable pair of diagrams $\cD_{g,g'}:f^o_1\to f^o_2$
and $\cD_{h,h'}:f^o_2\to f^o_3$.  However, by proposition
\ref{prop_square-algebra} we have :
$$
\tau^{g,g'}\boxminus\tau^{h,h'}=
(\one_{g_*\circ\sigma_{f_2}}\boxminus\one_{h_*\circ\sigma_{f_3}})\boxvert
(\Upsilon(\cD_{g,g'})\boxminus\Upsilon(\cD_{h,h'}))\boxvert
(\omega^{g'}\boxminus\omega^{h'}).
$$
Clearly $\one_{g_*\circ\sigma_{f_2}}\boxminus\one_{h_*\circ\sigma_{f_3}}=
\one_{(g\circ h)_*\circ\sigma_{f_3}}$, and we have
$\omega^{g'}\boxminus\omega^{h'}=\omega^{g'\circ h'}$, by the
uniqueness properties of the isomorphisms $\omega^\bullet$.
Lastly, $\Upsilon(\cD_{g,g'})\boxminus\Upsilon(\cD_{h,h'})=
\Upsilon(\cD_{g,g'}\boxvert\cD_{h,h'})$, by proposition
\ref{prop_composition-of-sqlinks}(i), whence the claim.
\end{pfclaim}

From claim \ref{cl_patience}(ii) and remark
\ref{rem_functors-from-lax}(i) we deduce a functor
$\cFib(\Sigma^X_\bullet):\cC\to\Set$.

Next, notice that, for every object $(S'\xrightarrow{f}S)^o$
of $\sMorph(\Topos^o)$, the objects of the fibre category
$\sMorph(\pi)^{-1}(f)$ are the triples $(\zeta,\zeta',\beta)$,
where $\zeta$ is an $S$-point, $\zeta'$ is an $S'$-point, and
$\beta:\zeta\isom f(\zeta')$ is an isomorphism of $S$-points. The
morphisms $(\zeta_1,\zeta'_1,\beta_1)\to(\zeta_2,\zeta'_2,\beta_2)$
are the pairs $(\alpha,\alpha')$ where $\alpha:\zeta_1\isom\zeta_2$
and $\alpha':\zeta'_1\isom\zeta'_2$ are isomorphisms of $S$-points
and respectively $S'$-points, such that
\set\begin{equation}\label{eq_such-that}
(f_**\alpha')\odot\beta_1=\beta_2\odot\alpha.
\end{equation}
If $(\gamma,\gamma'):(\zeta_2,\zeta'_2,\beta_2)\to
(\zeta_3,\zeta'_3,\beta_3)$ is another such morphism, we have
$$
(\gamma,\gamma')\circ(\alpha,\alpha')=
(\gamma\odot\alpha,\gamma'\odot\alpha')
\qquad
\text{in $\sMorph(\pi)^{-1}(f)$}.
$$
The (split) cleavage $\sc$ of $\sMorph(\cFib(\Pt^\times))$ attaches
to the diagram $\cD_{g,g'}$, the functor
$$
\sc_{g,g'}:
\sMorph(\pi)^{-1}(f_2)\to\sMorph(\pi)^{-1}(f_1)
\qquad
(\zeta,\zeta',\beta)\mapsto(g(\zeta),g'(\zeta'),g(\beta))
$$
that assigns to every morphism $(\alpha,\alpha')$ of
$\sMorph(\pi)^{-1}(f_2)$ the morphism $(g(\alpha),g'(\alpha'))$.
With this notation, $p$ is the (cartesian) functor associated
with the strict pseudo-natural transformation
$$
t:\sc\Rightarrow\Pt^\times\circ\ss^o
$$
that assigns to every morphims of topoi $f:S'\to S$
the functor
$$
\sMorph(\pi)^{-1}(f)\to\Pt(S)
\qquad
(\zeta,\zeta',\beta)\mapsto\zeta
\qquad
(\alpha,\alpha')\mapsto\alpha.
$$
Now, to every morphism of topoi $f:S'\to S$, and every
$(\zeta,\zeta',\beta)\in\Ob(\sMorph(\pi)^{-1}(f))$, let
us attach the map
$b_{f,(\zeta,\zeta',\beta)}:\Gamma(S',f^*\zeta_*X)\to\zeta^*\zeta_*X$
defined as the composition :
$$
\Gamma(S',f^*\zeta_*X)
\xrightarrow{\Gamma(S',(\eta^{\zeta'}*f^*\zeta_*)_X)}
\Gamma(S',\zeta'_*\zeta'^*f^*\zeta_*X)
\xrightarrow{\omega^{\zeta'}_{\zeta'^*f^*\zeta_*X}}
\zeta'^*f^*\zeta_*X\xrightarrow{(\beta^\dagger*\zeta_*)_X}
\zeta^*\zeta_*X
$$
where $\eta^{\zeta'}:\one_{S'}\Rightarrow\zeta'_*\zeta'^*$
is the unit of adjunction, and
$\omega^{\zeta'}:\Gamma_{S'}\circ\zeta'\isom\one_\Set$ is
the unique isomorphism of morphisms of topoi provided
by proposition \ref{prop_glob-sections}(ii).

\begin{claim}\label{cl_sharecropping}
(i)\ \
For every morphism of topoi $f:S'\to S$, the rule :
$(\zeta,\zeta',\beta)\mapsto b_{f,(\zeta,\zeta',\beta)}$
defines a natural transformation
$$
b_f:\Sigma^X_f\circ t_f\Rightarrow\alpha^X_S\circ t_f.
$$

(ii)\ \
The rule : $f\mapsto b_f$ defines a modification
$$
b_\bullet:(\alpha^X*\ss^o)\odot{}^ot\leadsto\Sigma^X_\bullet\odot{}^ot.
$$
\end{claim}
\begin{pfclaim}(i): Let $(\alpha,\alpha'):(\zeta_1,\zeta'_1,\beta_1)
\to(\zeta_2,\zeta'_2,\beta_2)$ be a morphism in
$\sMorph(\pi)^{-1}(f)$, and set
$U:=(\zeta'_{2*}*\alpha'^\dagger*f^*\zeta_{2*})^{-1}_X\circ
(\zeta'_{2*}*\alpha'*f^*\zeta_{2*})_X$. We notice that :
$$
\begin{aligned}
(\eta^{\zeta'_2}*f^*\zeta_{2*})_X\circ(f^*(\alpha_X)&\,=
U\circ(\eta^{\zeta'_1}*f^*\zeta_{2*})_X\circ(f^*(\alpha_X) \\
&\,=U\circ(\zeta'_{1*}\zeta'^*_1f^**\alpha)_X\circ
(\eta^{\zeta'_1}*f^*\zeta_{1*})_X
\end{aligned}
$$
where the first equality holds by remark
\ref{rem_adjoint-transf}(iii), and the second one follows
from \eqref{eq_Godement-functor}. We are then reduced to
checking the commutativity of the diagram :
$$
\xymatrix@C+50pt{ \Gamma(S',\zeta'_{1*}\zeta'^*_1f^*\zeta_{1*}X)
\ar[r]^-{\omega^{\zeta'_1}_{\zeta'^*_1f^*\zeta_{1*}X}}
\ar[d]_{\Gamma(S',\zeta'_{1*}\zeta'^*_1f^*\alpha_X)} &
\zeta'^*_1f^*\zeta_{1*}X \ar[r]^-{(\beta^\dagger_1*\zeta_{1*})_X}
\ar[d]^{\zeta'^*_1f^*\alpha_X} &
\zeta^*_1\zeta_{1*}X \ar[d]^{\zeta^*_1\alpha_X} \\
\Gamma(S',\zeta'_{1*}\zeta'^*_1f^*\zeta_{2*}X)
\ar[r]^-{\omega^{\zeta'_1}_{\zeta'^*_1f^*\zeta_{2*}X}}
\ar[d]_{\Gamma(S',(\alpha'*\zeta'^*_1f^*\zeta_{2*})_X)} &
\zeta'^*_1f^*\zeta_{2*}X \ar[r]^-{(\beta^\dagger_1*\zeta_{2*})_X}
\ddouble & \zeta^*_1\zeta_{2*}X \ddouble \\
\Gamma(S',\zeta'_{2*}\zeta'^*_1f^*\zeta_{2*}X)
\ar[d]_{\Gamma(S',(\zeta'_{2*}*\alpha'^\dagger*f^*\zeta_{2*})^{-1}_X)}
\ar[r]^-{\omega^{\zeta'_2}_{\zeta'^*_1f^*\zeta_{2*}X}} &
\ar[d]^{(\alpha'^\dagger*f^*\zeta_{2*})^{-1}_X}
\zeta'^*_1f^*\zeta_{2*}X \ar[r]^-{(\beta^\dagger_1*\zeta_{2*})_X} &
\zeta^*_1\zeta_{2*}X \ar[d]^{(\alpha^\dagger*\zeta_{2*})^{-1}_X} \\
\Gamma(S',\zeta'_{2*}\zeta'^*_2f^*\zeta_{2*}X)
\ar[r]^-{\omega^{\zeta'_2}_{\zeta'^*_2f^*\zeta_{2*}X}} &
\zeta'^*_2f^*\zeta_{2*}X \ar[r]^-{(\beta^\dagger_2*\zeta_{2*})_X} &
\zeta^*_2\zeta_{2*}X. 
}$$
However, the commutativity of the top two squares and of
the bottom square on the left is clear; the commutativity
of the bottom square on the right follows from remark
\ref{rem_adjoint-transf}(iv,v) and \eqref{eq_such-that}.
Lastly, the commutativity of the left square on the central
row follows immediately from the uniqueness properties of
the isomorphisms $\omega^{\zeta'_1}$ and $\omega^{\zeta'_2}$
(proposition \ref{prop_glob-sections}(ii)), whence the
contention.

(ii): Consider a diagram $\cD_{g,g'}$ as in the foregoing, and
an object $(\zeta,\zeta',\beta)$ of $\sMorph(\pi)^{-1}(f_2)$.
We have to check the compatibility condition :
$$
b_{f_2,(\zeta,\zeta',\beta)}\circ\tau^{g,g'}_\zeta=
\tau^X_{g,\zeta}\circ b_{f_1,(g(\zeta),g(\zeta'),g(\beta))}.
$$
To this aim, notice that
$\omega^{g'(\zeta')}=\omega^{\zeta'}\odot(\omega^{g'}*\zeta'_*)$,
due to proposition \ref{prop_glob-sections}(ii); then it
suffices to show the commutativity of the diagram :
$$
\xymatrix@C-16pt@R+10pt{ \Gamma(S'_1,f_1^*g_*\zeta_*X)
\ar[rrr]^-{\Gamma(S'_1,\Upsilon(\cD_{g,g'})*\zeta_*)_X}
\ar[d]|{\Gamma(S'_1,(\eta^{g'(\zeta')}*f_1^*g_*\zeta_*)_X)} & & &
\Gamma(S'_1,g'_*f_2^*\zeta_*X)
\ar[rr]^-{\omega^{g'}_{f_2^*\zeta_*X}} & & \Gamma(S'_2,f_2^*\zeta_*X)
\ar[dd]|{\Gamma(S'_2,(\eta^{\zeta'}*f^*_2\zeta_*)_X)} \\
\Gamma(S'_1,g'(\zeta')_*g'(\zeta')^*f_1^*g_*\zeta_*X)
\ar[d]|{\omega^{g'}_{\zeta'_*\zeta'^*g'^*f_1^*g_*\zeta_*X}} \\
\Gamma(S'_2,\zeta'_*g'(\zeta')^*f_1^*g_*\zeta_*X)
\ar[d]|{\omega^{\zeta'}_{\zeta'^*g'^*f^*_1g_*\zeta_*X}}
\ar[rrrrr]^-{\Gamma(S'_2,(\zeta'^*f_2^**\eps^g*\zeta_*)_X)}
& & & & & \Gamma(S'_2,\zeta'_*\zeta'^*f_2^*\zeta_*X)
\ar[d]|{\omega^{\zeta'}_{\zeta'^*f^*_2\zeta_*X}} \\
\zeta'^*g'^*f_1^*g_*\zeta_*X \rdouble
\ar[d]|{(g(\beta)^\dagger*g(\zeta)_*)_X} & \zeta'^*f_2^*g^*g_*\zeta_*X
\ar[rrrr]^-{(\zeta'^*f_2^**\eps^g*\zeta_*)_X} & & & &
\zeta'^*f^*_2\zeta_*X \ar[d]|{(\beta^\dagger*\zeta_*)_X} \\
g(\zeta)^*g(\zeta)_*X \ar[rrrrr]^-{(\zeta^**\eps^g*\zeta_*)_X}
& & & & & \zeta^*\zeta_*X.
}$$
However, the commutativity of the bottom square follows from
\eqref{eq_Godement-functor} and remark \ref{rem_adjoint-transf}(vi),
and that of the central square is clear, by the naturality of
$\omega^{\zeta'}$. In order to show the commutativity of the
top square, set $A:=\zeta_*X$, and recall that
$$
\Upsilon(\cD_{g,g'})=(g'_*f^*_2*\eps^g)\odot(\eta^{g'}*f^*_1g_*)
\qquad\text{and}\qquad
\eta^{g'(\zeta')}=(g'_**\eta^{\zeta'}*g'^*)\odot\eta^{g'}
$$
(proposition \ref{prop_opp-links-and-base-ch}). We are then
reduced to checking the commutativity of the diagram :
$$
\xymatrix@R+15pt{ & \Gamma(S'_1,g'_*g'^*f_1^*g_*A)
\ar[ld]_-{\Gamma(S'_1,(g'_**\eta^{\zeta'}*g'^*f^*_1)_A)\ \ \ \ \ \ }
\ar[d]|{\omega^{g'}_{g'^*f^*_1g_*A}}
\ar[rrr]^-{\Gamma(S'_1,(g'_*f^*_2*\eps^g)_A)} & & &
\Gamma(S'_1,g'_*f^*_2A) \ar[d]|{\omega^{g'}_{f^*_2A}} \\
\Gamma(S'_1,g'(\zeta')_*g'(\zeta')^*f^*_1g_*A)
\ar[rd]_-{\omega^{g'}_{\zeta'_*g'(\zeta')^*f^*_1g_*A}} &
\Gamma(S'_2,g'^*f^*_1g_*A)
\ar[d]|{\Gamma(S'_2,(\eta^{\zeta'}*g'^*f^*_1g_*)_A)}
\ar[rrr]^-{\Gamma(S'_2,(f^*_2*\eps^g*)_A)} & & &
\Gamma(S'_2,f^*_2A) \ar[d]|{\Gamma(S'_2,(\eta^{\zeta'}*f^*_2)_A)} \\
& \Gamma(S'_2,\zeta'_*g'(\zeta')^*f^*_1g_*A)
\ar[rrr]^-{\Gamma(S'_2,(\zeta'_*f^*_2*\eps^g)_A)} & & &
\Gamma(S'_2,\zeta'_*\zeta'^*f^*_2A).
}$$
However, it is easily seen that the two square subdiagrams
and the triangular subdiagram commute (details left to the
reader), whence the contention.
\end{pfclaim}

From claim \ref{cl_sharecropping}(ii) and remark
\ref{rem_functors-from-lax}(ii) we get the natural
transformation
$$
\cFib(b_\bullet):
\cFib(\Sigma^X_\bullet)\circ p\Rightarrow\Xi^X\circ\ss'\circ p.
$$
Lastly, let $q:\cL(\xi)\to\cFib(\Pt^\times)$ be the natural
projection, $\ss_{\cL(\xi)}:\sMorph(\cL(\xi))\to\cL(\xi)$
the source functor, and set
$\sG:=\cFib(\Sigma^X_\bullet)\circ p\circ\sMorph(q)$; we
deduce the natural transformation
$$
\beta:\cFib(b_\bullet)*\sMorph(q):
\sG\Rightarrow\Xi^{X,\xi}\circ\ss_{\cL(\xi)}.
$$
By theorem \ref{th_Kan-ext}, the latter induces a natural
transformation
\set\begin{equation}\label{eq_apply-Kan}
\int^{\ss_{\cL(\xi)}}\sG\Rightarrow\Xi^{X,\xi}.
\end{equation}
Notice that, under the natural identification
$\sMorph(\cL(\xi))\isom\sMorph(\bLift(\xi))^o$,
the functor $\ss_{\cL(\xi)}$ correspond to the functor
$\st^o_{\bLift(\xi)}$, where
$\st_{\bLift(\xi)}:\sMorph(\bLift(\xi))\to\bLift(\xi)$
is the target functor. Moreover, $\st_{\bLift(\xi)}$ is a
fibration, by virtue of remark \ref{rem_Nbd-as-Fib}(ii)
and example \ref{ex_fibred-cats}(iii). Furthermore, for
every $(U,\xi_U,\omega_U)\in\Ob(\bLift(\xi))$, the fibre
category $\st^{-1}(U,\xi_U,\omega_U)$ is
$$
\bLift(\xi)/(U,\xi_U,\omega_U)\isom\bLift(\xi_U).
$$
Hence, the restriction of $\beta$ to the fibre category
$\st^{-1}(U,\xi_U,\omega_U)$ is a cocone
\set\begin{equation}\label{eq_messy-cocone}
\xymatrix@C+30pt{
\bLift(\xi_U)^o \rtwocell^{\sG}_{c_{\xi_U^*\xi_{U*}X}}{} & \Set.
}\end{equation}
By a direct inspection of the construction, we see that
\eqref{eq_messy-cocone} is the cocone $\tau_{\xi,\xi_{U*}X}*N^o_{\xi_U}$,
where $N_{\xi_U}:\bLift(\xi_U)\isom\bNbd(\xi_U)$ is the
equivalence \eqref{eq_same-neighbor}, and $\tau_{\xi,\xi_{U*}X}$
is the universal cocone of remark \ref{rem_Nbd-as-Fib}(iv).
Hence, \eqref{eq_messy-cocone} is a universal cocone, for
every $(U,\xi_U,\omega_U)\in\Ob(\bLift(\xi))$; taking
into account example \ref{ex_Fubini-fibred}(iv), we deduce
that \eqref{eq_apply-Kan} is an isomorphism of functors, and
we get a natural isomorphism :
$$
\colim_{\sMorph(\bLift(\xi))^o}\sG\isom\colim_{\bLift(\xi)^o}\Xi^{X,\xi}.
$$
Now, notice that the functor
\set\begin{equation}\label{eq_reveal}
\bLift(\xi)^o\to\sMorph(\bLift(\xi))^o
\qquad
(U,\xi_U,\omega_U)\mapsto\one_{(U,\xi_U,\omega_U)}
\end{equation}
is cofinal. A simple inspection reveals that the composition
of $\sG$ with \eqref{eq_reveal} is the constant functor with
value $X$, whence the contention.
\end{proof}

\sset\subsubsection{}\label{subsec_neighbors-site}
Suppose now that $T=C^\sim$ for some small site $C:=(\cC,J)$.
For a point $\xi$ of $T$, we may define another category
$\bNbd(\xi,C)$, whose objects are the pairs $(U,a)$ where
$U\in\Ob(\cC)$ and $a\in(h^a_U)_\xi$; the morphisms
$(U,a)\to(U',a')$ are the morphisms $f:U\to U'$ in $\cC$
such that $(h^a_f)_\xi(a)=a'$. (Notation of remark
\ref{rem_rep-and-sheafify}(iii).)
The rule $(U,a)\mapsto(h^a_U,a)$ defines a functor
\set\begin{equation}\label{eq_neighbors}
\bNbd(\xi,C)\to\bNbd(\xi).
\end{equation}

\begin{proposition}\label{prop_neighbors}
In the situation of \eqref{subsec_neighbors-site}, we have :
\begin{enumerate}
\item
The category $\bNbd(\xi,C)$ is cofiltered.
\item
The functor \eqref{eq_neighbors} is cofinal.
\end{enumerate}
\end{proposition}
\begin{proof} To begin with, since $\xi^*$ commutes with all
colimits, remark \ref{rem_rep-and-sheafify}(iii) implies easily
that, for every object $(U,a)$ of $\bNbd(\xi)$ there exist an
object $(V,b)$ of $\bNbd(\xi,C)$ and a morphism $(h^a_V,b)\to(U,a)$
in $\bNbd(\xi)$. Hence, it suffices to show (i).

Thus, let $(U_1,a_1)$ and $(U_2,a_2)$ be two objects of
$\bNbd(\xi,C)$; since $\bNbd(\xi)$ is cofiltered, we
may find an object $(F,b)$ of $\bNbd(\xi)$ and morphisms
$\phi_i:(F,b)\to(h^a_{U_i},a_i)$ (for $i=1,2$). By the foregoing,
we may also assume that $(F,b)=(h^a_V,b)$ for some object
$(V,b)$ of $\bNbd(\xi,C)$, in which case $\phi_i\in h^a_{U_i}(V)$
for $i=1,2$. By remark \ref{rem_rep-and-sheafify}(iv), we may
find a sieve $\cS$ covering $V$ such that
$\phi_i\in\Hom_{\cC^\wedge}(h_\cS,h^\sep_{U_i})$ for both
$i=1,2$.

On the other hand, \eqref{eq_colim-sieve} and proposition
\ref{prop_cover-is-iso} yield a natural isomorphism :
$$
\colim_\cS h^a\circ\ss\isom h^a_V.
$$
Therefore, since $\xi^*$ commutes with all colimits, we may
find $(f:S\to V)\in\Ob(\cS)$ and $c\in(h^a_S)_\xi$ such that
$h_f^a:(h_S^a,c)\to(h^a_V,b)$ is a morphism of neighborhoods
of $\xi$.

For $i=1,2$, denote by $\bar\phi_{S,i}\in h^\sep_U(S)$ the
image of $\phi_i$ (under the map induced by the natural
morphism $h_S\to h_\cS$ coming from \eqref{eq_colim-sieve}).
Pick any $\phi_{S,i}\in\Hom_\cC(S,V)$ in the preimage of
$\bar\phi_{S,i}$; then $\phi_{S,i}$ defines a morphism
$(S,c)\to(U_i,a_i)$ in $\bNbd(\xi,C)$.

Next, suppose that $\phi_1,\phi_2:(U,a)\to(U',a')$
are two morphisms in $\bNbd(\xi,C)$; arguing as in the foregoing,
we may find an object $(V,b)$ of $\bNbd(\xi,C)$, and
$\psi\in h_U^\sep(V)=\Hom_{\cC^\wedge}(h^\sep_V,h_U^\sep)$
whose image in $h^a_U(V)$ yields a morphism
$\psi^a:(h^a_V,b)\to(h^a_U,a)$ in $\bNbd(\xi)$,
and such that $\phi_1^\sep\circ\psi=\phi_2^\sep\circ\psi$
in $h^\sep_{U'}(V)$.
We may then find a covering subobject $i:R\to h_V$, such that
$\phi_1\circ(\psi\circ i)=\phi_2\circ(\psi\circ i)$ in
$\Hom_{\cC^\wedge}(R,h_U')$. Again, by combining
\eqref{eq_colim-sieve} and proposition \ref{prop_cover-is-iso},
we deduce that there exists a morphism
$\beta:(V',b')\to(V,b)$ in $\bNbd(\xi,C)$, such that
$\phi_1\circ(\psi\circ\beta)=\phi_2\circ(\psi\circ\beta)$.
This completes the proof of (i).
\end{proof}

As a corollary of proposition \ref{prop_neighbors} and
of remark \eqref{rem_Nbd-as-Fib}(iv), we deduce, for every
sheaf $F$ on $C$, a natural isomorphism :
$$
\colim_{\bNbd(\xi,C)^o}F\circ\iota^o_C\isom F_\xi
$$
where $\iota_C:\bNbd(\xi,C)\to\cC$ is the functor such that
$(U,a)\mapsto U$ for every object $(U,a)$.

\subsection{Algebra on a topos}\label{sec_tensor-on-topoi}
Let $T$ be any topos, and endow $T$ with the structure of
tensor category as explained in example \ref{ex_stupid-tensor}
(so, the tensor functor is given by fixed choices of products
for every pairs of objects of $T$, and any final object $1_T$
can be taken for unit object of $(T,\otimes)$). We notice that
$(T,\otimes)$ admits an internal $\Hom$ functor (see remark
\ref{rem_Hom-how-to}(ii)). Indeed, let $X$ and $X'$ be any
two objects of $T$. It is easily seen that the presheaf on $T$ :
$$
U\mapsto\Hom_{T\!/U}(X'_{|U},X_{|U})=\Hom_T(X'\!\times\!U,X)
$$
is actually a sheaf on $(T,C_T)$ (notation of example
\ref{ex_localization-topos}(i)), so it is an object of $T$,
denoted :
$$
\cHom_T(X',X).
$$
The functor :
$$
T\to T \quad : \quad X\mapsto\cHom_T(X',X)
$$
is right adjoint to the functor $T\to T$ : $Y\mapsto Y\times X'$,
so it is an internal $\Hom$ functor for $X'$.

If $f:T\to S$ is a morphism of topoi, and $Y\in\Ob(S)$, we have
a natural isomorphism in $S$ :
\set\begin{equation}\label{eq_projection_form} \cHom_S(Y,f_*X)\isom
f_*\cHom_T(f^*Y,X)
\end{equation}
which, on every $U\in\Ob(S)$, induces the natural bijection :
$$
\Hom_S(Y\!\times\!U,f_*X)\isom\Hom_T(f^*Y\!\times\!f^*U,X)
$$
given by the adjunction $(f^*,f_*)$. By general nonsense, from
\eqref{eq_projection_form} we derive a natural morphism in $S$ :
$$
f_*\cHom_T(X',X)\to\cHom_S(f_*X',f_*X)
$$
and in $T$ :
$$
f^*\cHom_S(Y',Y)\xrightarrow{ \theta_f }\cHom_T(f^*Y',f^*Y) \qquad
\text{for any $Y,Y'\in\Ob(S)$.}
$$
Moreover, if $g:U\to T$ is another morphism of topoi, the diagram :
\set\begin{equation}\label{eq_obvious-but-cumber}
{\diagram
g^*f^*\cHom_S(Y',Y) \ar[rr]^{g^*\theta_f}
\ar[rd]_{\theta_{f\circ g}} & & g^*\cHom_T(f^*Y',f^*Y)
\ar[ld]^{\theta_g} \\
& \cHom_U(g^*f^*Y',g^*f^*Y) \enddiagram}
\end{equation}
commutes, up to a natural isomorphism.

\sset\subsubsection{}
Let $A$ be any object of $T$, and $(X,\mu_X)$ a left
$A$-module for the tensor category structure on $T$ as in
\eqref{sec_tensor-on-topoi}; for every object $U$ of $T$ we
obtain a left $A_{|U}$-module (on the topos $T/U$ : see example
\ref{ex_localization-topos}(i)), by the rule :
$$
(X,\mu_X)_{|U}:=(X_{|U},\mu_X\times\one_U).
$$
If $(X',\mu_{X'})$ is another left $A$-module, it is easily seen
that the presheaf on $T$ :
$$
U\mapsto\Hom_{A_{|U}\Mod_l}((X,\mu_X)_{|U},(X',\mu_{X'})_{|U})
$$
is actually a sheaf for the canonical topology, so it is an object
of $T$, denoted :
$$
\cHom_{A_l}((X,\mu_X),(X',\mu_{X'}))
$$
(or just $\cHom_{A_l}(X,X')$, if the notation is not ambiguous). The
same considerations can be repeated for the sets of morphisms of
right $B$-modules, and of $(A,B)$-bimodules, so one gets objects
$\cHom_{B_r}(X,X')$ and $\cHom_{(A,B)}(X,X')$. By a simple inspection,
we see that these objects are naturally isomorphic to the objects
denoted in the same way in \eqref{subsec_Hom-as-equal}, so the
notation is not in conflict with {\em loc.cit.}; it also follows
that $\cHom_{A_l}(X,X')$ is the equalizer of two morphisms in $T$ :
$$
\xymatrix{\cHom_T(X,X') \ar@<.5ex>[r] \ar@<-.5ex>[r]
& \cHom_T(A\times X,X')}.
$$
In the same vein, let $A,B,C\in\Ob(T)$ be any three objects,
$S$ an $(A,B)$-bimodule, $S'$ a $(C,B)$-bimodule, and $S''$
a $(C,A)$-bimodule. Then the $(C,B)$-bimodule
$\cHom_{B_r}(S,S')$ and the $(C,B)$-bimodule $S''\otimes_AS$
(see \eqref{subsec_adjunct-scoppia}) are the sheaves on
$(T,C_T)$ associated with the presheaves given by the rules :
$U\mapsto\Hom_{B_{|U,r}}(S_{|U},S'_{|U})$, and respectively :
$U\mapsto S'(U)\otimes_{M(U)}S(U)$ for every object $U$ of
$T$. Furthermore, the general theory of monoids, their
modules and their tensor products, developed in section
\ref{sec_tensors-ab} is available in the present situation,
so we have a well defined notion of $T$-monoid
(see example \ref{ex_stupid-prior}(ii) and remark
\ref{rem_commutes-forgets}). Via the equivalence of theorem
\eqref{th_canon-topos}, a $T$-monoid $\underline M$ is
also the same as a sheaf of monoids $M$ on the site
$(T,C_T)$, and a left (resp. right, resp. bi-)
$\underline M$-module is the same as the datum of a
sheaf $S$ in $(T,C_T)$, such that $S(U)$ is a left
(resp. right, resp. bi-) $M(U)$-module, for every
object $U$ of $T$.

\sset\subsubsection{}\label{subsec_up_and-down-mods}
Let $f:T_1\to T_2$ be a morphism of topoi, $A_1$ an object
of $T_1$, and $(X,\mu_X)$ a left $A_1$-module.
Since $f_*$ is left exact, we have a natural isomorphism :
$f_*(A_1\times X)\isom f_*A_1\times f_*X$, so we obtain a left
$f_*A_1$-module :
$$
f_*(X,\mu_X):=(f_*X,f_*\mu_X)
$$
which we denote just $f_*X$, unless the notation is ambiguous.
Likewise, since $f^*$ is left exact, from any object $A_2$ of
$T_2$, and any left $A_2$-module $(Y,\mu_Y)$, we obtain a left
$f^*A_2$-module :
$$
f^*(Y,\mu_Y):=(f^*Y,f^*\mu_Y).
$$
The same considerations apply of course, also to right modules and
to bimodules. Furthermore, let $A_i,B_i,C_i$ be three objects of
$T_i$ (for $i=1,2$); since $f_*$ is left exact, for any
$(A_1,B_1)$-bimodule $X$ and any $(C_1,A_1)$-bimodule $X'$ we have
a natural morphism of $(f_*C_1,f_*B_1)$-bimodules :
$$
f_*X'\otimes_{f_*A_1}f_*X\to f_*(X'\otimes_{A_1} X)
$$
and since $f^*$ is exact, for any $(A_2,B_2)$-bimodule
$Y$ and any $(C_2,A_2)$-bimodule $Y'$ we have a natural
isomorphism :
$$
f^*Y'\otimes_{f^*A_2}f^*Y\isom f^*(Y'\otimes_{A_2}Y)
$$
of $(f^*C_2,f^*B_2)$-bimodules.

\sset\subsubsection{}\label{eq_presheaves-mods}
The constructions of the previous paragraphs also apply to
presheaves on $T$ : this can be seen, {\em e.g.} as follows.
Pick a universe $\sV$ such that $T$ is $\sV$-small; then
$T^\wedge_\sV$ is a $\sV$-topos. Hence, if $A$ is any
$\sV$-presheaf on $T$, and $X,X'\in\Ob(T^\wedge_\sV)$ two
left $A$-modules, we may construct $\cHom_{A_l}(X,X')$ as
an object in $T^\wedge_\sV$. Now, if $A,X,X'$ lie in the
full subcategory $T^\wedge_\sU$ of $T^\wedge_\sV$, it is
easily seen that also $\cHom_{A_l}(X,X')$ lies in $T^\wedge_\sU$.
Likewise, if $A,B,C$ are two $\sU$-presheaves on $T$, we
may define $X'\otimes_AX$ in $T^\wedge_\sU$, for any
$(A,B)$-bimodule $X$ and $(C,A)$-bimodule $X'$, and this
tensor product will still be left adjoint to the
$\cHom$-functor for presheaves on $T$.

Moreover, we have a natural morphism of topoi
$i_\sV:(T,C_T)^\sim_\sV\to T^\wedge_\sV$, given by the
forgetful functor and its left adjoint $i_\sV^*$, which is
given by the rule : $F\mapsto F^a$ : see example
\ref{ex_another-basic}(i). The restriction of $i_\sV^*$ to
the full subcategory $T^\wedge_\sU$ is isomorphic to a functor
that factors through the Yoneda embedding $T\to(T,C_T)^\sim_\sV$,
therefore the discussion of \eqref{subsec_up_and-down-mods}
specializes to show that, for every $\sU$-presheaf $A$
on $T$, and every $A$-module $X\in\Ob(T^\wedge_\sU)$,
the object $X^a\in\Ob(T)$ is naturally an $A^a$-module.
Also, for any $(A,B)$-bimodule $X$ and $(C,A)$-bimodule $X'$,
such that $A,B,C,X,X'$ are $\sU$-small, we have a natural
isomorphism of $(C^a,B^a)$-bimodules :
$$
X^{\prime a}\otimes_{A^a}X^a\isom(X'\otimes_AX)^a.
$$
The following definition gathers some further notions
-- specific to monoids over a topos -- which shall be
used in this work.

\begin{definition}\label{def_coh-idea-log}
Let $T$ be a topos, $\underline M$ a $T$-monoid, $S$ a left
(resp. right, resp. bi-) $\underline M$-module.
\begin{enumerate}
\item
$S$ is said to be {\em of finite type}, if there exist a
covering family $(U_\lambda\to 1_T~|~\lambda\in\Lambda)$ of the
final object of $T$, and for every $\lambda\in\Lambda$ an integer
$n_\lambda\in\N$ and an epimorphism of left (resp. right, resp. bi-)
$\underline M_{|U_\lambda}$-modules :
$\underline M^{\oplus n_\lambda}_{|U_\lambda}\to S_{|U_\lambda}$.
\item
$S$ is {\em finitely presented}, if there exist a covering
family $(U_\lambda\to 1_T~|~\lambda\in\Lambda)$ of the
final object of $T$, and for every $\lambda\in\Lambda$ integers
$m_\lambda,n_\lambda\in\N$ and morphisms
$f_\lambda,g_\lambda:
\underline M^{\oplus m_\lambda}_{|U_\lambda}\to
\underline M^{\oplus n_\lambda}_{|U_\lambda}$ whose coequalizer
-- in the category of left (resp. right, resp. bi-)
$\underline M_{|U_\lambda}$-modules -- is isomorphic to
$S_{|U_\lambda}$.
\item
$S$ is said to be {\em coherent}, if it is of finite type,
and for every object $U$ in $T$, every submodule of finite type
of $S_{|U}$ is finitely presented.
\item
$S$ is said to be {\em invertible}, if there exist a covering
family $(U_\lambda\to 1_T~|~\lambda\in\Lambda)$, and for every
$\lambda\in\Lambda$, an isomorphism
$\underline M_{|U_\lambda}\isom S_{|U_\lambda}$ of left (resp.
right, resp. bi-) $\underline M_{|U}$-modules. (Thus, every
invertible module is finitely presented.)
\item
An ideal $I\subset\underline M$ is said to be {\em invertible},
(resp. {\em of finite type}, resp. {\em finitely presented},
resp. {\em coherent}) if it is such, when regarded as an
$\underline M$-bimodule.
\end{enumerate}
\end{definition}

\begin{example}\label{ex_see-also}
(i)\ \
Take again $T=\Set$. Then an $M$-module $S$ is of finite type
if and only if there exists a finite subset $\Sigma\subset S$,
such that $S=M\cdot\Sigma$, with obvious notation.
In this case, we say that $\Sigma$ is a {\em finite system of
generators\/} of $S$ (and we also say that $S$ is
{\em finitely generated}; likewise, an ideal of finite type
is also called finitely generated). We say that $S$ is
{\em cyclic}, if $S=M\cdot s$ for some $s\in S$.

(ii)\ \
If $S$ and $S'$ are two $M$-modules, the coproduct
$S\oplus S'$ is the disjoint union of $S$ and $S'$, with scalar
multiplication given by the disjoint union of the laws $\mu_S$
and $\mu_{S'}$. The product $S\times S'$ is the cartesian
product of the underlying sets, with scalar multiplication given by
the rule : $x\cdot(s,s'):=(x\cdot s,x\cdot s')$ for every $x\in M$,
$s,\in S$, $s'\in S'$.
\end{example}

For future use, let us also make the :

\begin{definition}\label{def_grad-monoids}
Let $T$ be any topos, $\underline P$ a $T$-monoid, and $(N,+,0)$
any monoid.
\begin{enumerate}
\item
We say that $\underline P$ is {\em $N$-graded}, if it admits
a morphism of monoids $\pi:\underline P\to N_T$, where $N_T$
is the constant sheaf of monoids arising from $N$ (the
coproduct of copies of the final object $1_T$ indexed by $N$).
For every $n\in N$ we let $\underline P_n:=\pi^{-1}(n_T)$, the
preimage of the global section corresponding to $n$. Then
$$
\underline P=\coprod_{n\in N}\underline P_n
$$
the coproduct of the objects $\underline P_n$, and the
multiplication law of $\underline P$ restricts to a map
$\underline P_n\times\underline P_m\to\underline P_{n+m}$, for
every $n,m\in N$. Especially, each $\underline P_n$ is a
$\underline P_0$-module, and $\underline P$ is also the direct
sum of the $\underline P_n$, in the category of
$\underline P_0$-modules. The morphism $\pi$ is called the
{\em grading\/} of $\underline P$.
\item
In the situation of (i), let $S$ be a left (resp. right, resp. bi-)
$\underline P$-module. We say that $S$ is {\em $N$-graded}, if it
admits a morphism of $\underline P$-modules $\pi_S:S\to N_T$, where
$N_T$ is regarded as a $\underline P$-bimodule via the grading
$\pi$ of $\underline P$. Then $S$ is the coproduct
$S=\coprod_{n\in N}S_n$, where $S_n:=\pi_S^{-1}(n_T)$,
and the scalar multiplication of $S$ restricts to morphisms
$\underline P_n\times S_m\to S_{n+m}$, for every $n,m\in N$.
The morphism $\pi_S$ is called the {\em grading\/} of $S$.
\item
A morphism $\underline P\to\underline Q$ of $N$-graded
$T$-monoids is a morphism of monoids that respects the
gradings, with obvious meaning. Likewise one defines morphisms
of $N$-graded $\underline P$-modules.
\end{enumerate}
\end{definition}

\begin{example}
Take $T=\Set$, and let $M$ be any commutative monoid. Then we
claim that the only invertible object in the tensor category
$M\Mod_l$ is $M$; {\em i.e.} if $S$ and $S'$ are any two
$(M,M)$-bimodules, then $S\otimes_MS'\simeq M$ if and only
if $S$ and $S'$ are both isomorphic to $M$.

Indeed, let $\phi:S\otimes_MS'\isom M$ be an isomorphism,
and choose $s_0\in S$, $s'_0\in S'$ such that
$\phi(s_0\otimes s'_0)=1$. Consider the morphisms of left
$M$-modules :
$$
M\xrightarrow{\ \alpha\ }S\xrightarrow{\ \beta\ }M \qquad
M\xrightarrow{\ \alpha'\ }S'\xrightarrow{\ \beta'\ }M
$$
such that :
$$
\alpha(m)=m\cdot s_0 \quad \beta(s)=\phi(s\otimes s'_0) \quad
\alpha'(m)=m\cdot s'_0 \quad \beta'(s')=\phi(s_0\otimes s')
$$
for every $m\in M$, $s\in S$, $s'\in S'$; we notice that
$\beta\circ\alpha=\one_M=\beta'\circ\alpha'$. There follows
natural morphisms :
$$
S'\xrightarrow{\gamma}
M\otimes_MS'\xrightarrow{\alpha\otimes_MS'}S\otimes_MS'
\xrightarrow{\beta\otimes_MS'}M\otimes_MS'\xrightarrow{\gamma^{-1}}S'
$$
whose composition is the identity $\one_{S'}$. However, it is easily
seen that $\phi\circ(\alpha\otimes_MS')\circ\gamma=\beta'$ and
$\gamma^{-1}\circ(\beta\otimes_MS')\circ\phi^{-1}=\alpha'$, thus
$\alpha'\circ\beta'=\one_{S'}$, hence both $\alpha'$ and $\beta'$
are isomorphisms, and the same holds for $\alpha$ and $\beta$.
\end{example}

\begin{example}\label{ex_symm-pow}
Let $\underline M$ be a $T$-monoid, and $\cL$ a
$\underline M$-bimodule. For every $n\in\N$, let
$\cL^{\otimes n}:=
\cL\otimes_{\underline M}\cdots\otimes_{\underline M}\cL$,
the $n$-fold tensor power of $\cL$. The $\N$-graded
$\underline M$-bimodule
$$
\mathrm{Tens}^\bullet_{\underline M}\cL:=\coprod_{n\in\N}\cL^{\otimes n}
$$
is naturally a $\N$-graded $T$-monoid, with composition law
induced by the natural morphisms
$\cL^{\otimes n}\otimes_{\underline M}\cL^{\otimes m}\isom
\cL^{\otimes n+m}$, for every $n,m\in\N$. (Here we set
$\cL^{\otimes 0}:=\underline M$.) If $\cL$ is invertible,
$\mathrm{Tens}^\bullet_{\underline M}\cL$ is a commutative
$\N$-graded $T$-monoid, which we also denote
$\Sym^\bullet_{\underline M}\cL$.
\end{example}

\begin{remark}\label{rem_push-pull-mond}
(i)\ \ 
Let $f:T\to S$ be a morphism of topoi, $\underline M:=(M,\mu_M)$
a $T$-semigroup, and $\underline N:=(N,\mu_N)$ a $S$-semigroup.
Then clearly $f_*\underline M:=(f_*M,f_*\mu_M)$ is a $S$-semigroup,
and $f^*\underline N:=(f^*N,f^*\mu_N)$ is a $T$-semigroup.

(ii)\ \ 
Furthermore, if $1_M:1_T\to M$ (resp. $1_N:1_S\to N$)
is a unit for $M$ (resp. for $N$), then notice that $f_*1_T=1_S$
(resp. $f^*1_S=1_T$), since the final object is the empty product;
it follows that $f_*1_M$  (resp. $f^*1_S$) is a unit for
$f_*\underline M$ (resp. for $f^*\underline N$).

(iii)\ \ 
Obviously, if $\underline M$ (resp. $\underline N$) is
commutative, the same holds for $f_*\underline M$ (resp.
$f^*\underline N$).

(iv)\ \ 
If $X$ is a left $\underline M$-module, then $f_*X$
is a left $f_*\underline M$-module, and if $Y$ is a left
$\underline N$-module, then $f^*Y$ is a left $f^*\underline N$-module.
The same holds for right modules and bimodules.

(v)\ \ 
Moreover, let $\eps_M:f^*f_*\underline M\to\underline M$
(resp. $\eta_N:\underline N\to f_*f^*\underline N$) be the counit
(resp. unit) of adjunction. Then the counit (resp. unit) :
$$
\eps_X:f^*f_*X\to X_{(\eps_M)}
\qquad
\text{(resp. $\eta_Y:Y\to f_*f^*Y_{(\eta_N)}$)}
$$
is a morphism of $f^*f_*\underline M$-modules (resp. of
$\underline N$-modules) (notation of \eqref{subsec_tens-restr}).
(Details left to the reader.)

(vi)\ \
Let $\phi:f^*\underline N\to\underline M$ be a morphism of
$T$-monoids. Then the functor
$$
\underline N\Mod_l\to\underline M\Mod_l
\quad : \quad
Y\mapsto M\otimes_{f^*N}f^*Y
$$
is left adjoint to the functor :
$$
\underline M\Mod_l\to\underline N\Mod_l
\quad : \quad
X\mapsto f_*X_{(\eta_N)}.
$$
(And likewise for right modules and bimodules : details left
to the reader.)

(vii)\ \
The considerations of \eqref{eq_presheaves-mods} also
apply to monoids : we get that, for any presheaf of monoids
$\underline M:=(M,\mu_M,1_M)$ on $T$, the datum
$\underline M^a:=(M^a,\mu_M^a,1_M^a)$ is a $T$-monoid, and
we have a well defined functor :
$$
\underline M\Mod_l\to\underline M^a\Mod_l
\qquad
X\mapsto X^a.
$$
(And  as usual, the same applies to right modules and bimodules.)
\end{remark}

\sset\subsubsection{}\label{subsec_ext-mod-byempty}
Let $T$ be a topos, $U$ any object of $T$, and $\underline M$
a $T$-monoid. As a special case of remark
\ref{rem_push-pull-mond}(i), we have the $T\!/U$-monoid
$j_U^*\underline M=\underline M{}_{|U}$, and if we take
$\phi:=\one_{j_U^*\underline M}$ in remark
\ref{rem_push-pull-mond}(vi), we deduce that the functor
$$
j_U^*:\underline M\Mod_l\to\underline M{}_{|U}\Mod_l
\qquad
Y\mapsto Y_{|U}
$$
admits the right adjoint $j_{U*}$. Now, suppose that
$X\to U$ is any left $\underline M{}_{|U}$-module. The scalar
multiplication of $X$ is a $U$-morphism $\mu_X:M\times X\to X$
and $j_{U!}\mu_X$ is the same morphism, seen as a morphism
in $T$ (notation of example \ref{ex_localization-topos}(i)).
In other words, $j_{U!}$ induces a faithful functor on left
modules, also denoted :
$$
j_{U!}:\underline M{}_{|U}\Mod_l\to\underline M\Mod_l.
$$
It is easily seen that this functor is left adjoint to
the foregoing functor $j^*_U$. Especially, this functor
is right exact; it is not generally left exact, since
it does not preserve the final object (unless $U=1_T$).
However, it does commute with fibre products, and therefore
transforms monomorphisms into monomorphisms. All this holds
also for right modules and bimodules.

\sset\subsubsection{}\label{subsec_pointed-modules}
Let $T$ be any category as in example \ref{ex_stupid-prior}(i),
denote by $1_T$ a final object of $T$, and by $\underline M$
any $T$-monoid. A {\em pointed left $\underline M$-module\/} is
a datum
$$
(S,0_S)
$$
consisting of a left $\underline M$-module $S$ and a morphism of
$\underline M$-modules $0_S:1_T\to S$, where $0$ is the final object
of $M\Mod_l$.
Often we shall write $S$ instead of $(S,0_S)$, unless this may give
rise to ambiguities. As usual, a morphism $\phi:S\to T$ of pointed
modules is a morphism of $M$-modules, such that $0_T=\phi\circ 0_S$.
In other words, the resulting category is just $0/M\Mod_l$, and shall
be denoted $\underline M\Mod_{l\circ}$.

Likewise one may define the category $\underline M\Mod_{r\circ}$ of
right $\underline M$-modules, and $(\underline M,\underline N)\Mod_\circ$
of pointed bimodules, for given $T$-monoids $\underline M$ and
$\underline N$.

\begin{remark}\label{rem_apparent-reasons}
Let $T$ be a category as in remark \ref{rem_commutes-forgets}.

(i)\ \ 
For reasons that will become readily apparent, for many
purposes the categories of pointed modules are more useful
than the non-pointed variant of \eqref{subsec_sheaves-of-mons}.
In any case, we have a faithful functor :
\set\begin{equation}\label{eq_faithf-point}
\underline M\Mod_l\to\underline M\Mod_{l\circ}
\qquad
S\mapsto S_\circ:=(S\oplus 0,0_S)
\end{equation}
where $0_S:0\to S\oplus 0$ is the obvious inclusion map. Thus,
we may -- and often will, without further comment -- regard any
$\underline M$-module as a pointed module, in a natural way.
(The same can of course be repeated for right modules and bimodules.)

(ii)\ \ 
In turn, when dealing with pointed $M$-modules, things often work
out nicer if $\underline M$ itself is a {\em pointed $T$-monoid\/}.
The latter is the datum $(\underline M,0_M)$ of a $T$-monoid
$\underline M$ and a morphism of $\underline M$-modules
$0_M:0\to M$. A morphism of pointed $T$-monoids is of course just
a morphism $f:\underline M\to\underline M'$ of $T$-monoids, such
that  $f\circ 0_M=0_{M'}$. As customary, we shall often just write
$\underline M$ instead of $(\underline M,0_M)$, unless we wish
to stress that $\underline M$ is pointed.

(iii)\ \ 
Let $(\underline M,0_M)$ be a pointed $T$-monoid; a
{\em pointed left $(\underline M,0_M)$-module\/} is a pointed
left $\underline M$-module $S$, such that $0\cdot s=0$ for every
$s\in S$. A morphism of pointed left $(\underline M,0_M)$-modules
is just a morphism of pointed left $\underline M$-modules.
As usual, these gadgets form a category
$(\underline M,0_M)\Mod_{l\circ}$. Similarly we have
the right and bi-module variant of this definition.

(iv)\ \ 
The forgetful functor from the category of pointed $T$-monoids
to the category of $T$-monoids, admits a left adjoint :
$$
\underline M\mapsto(\underline M{}_\circ,0_{M_\circ}).
$$
Namely, $M_\circ$ is the $\underline M$-module $M\oplus 0$, the
zero map $0_{M_\circ}:0\to M\oplus 0$ is the obvious inclusion, and
the scalar multiplication $M\times M_\circ\to M_\circ$ is extended
to a multiplication law $\mu:M_\circ\times M_\circ\to M_\circ$ in the
unique way for which $(\underline M{}_\circ,\mu,0_{M_0})$ is a pointed
monoid. The unit of adjunction $\underline M\to\underline M{}_\circ$
is the obvious inclusion map.

(v)\ \ 
If $\underline M$ is a (non-pointed) monoid, the restriction of scalars
$$
(\underline M{}_\circ,0_{M_\circ})\Mod_{l\circ}\to
\underline M\Mod_{l\circ}
$$
is an isomorphism of categories. Namely, any pointed left
$\underline M$-module $S$ is naturally a pointed left
$\underline M{}_\circ$-module : the given scalar multiplication
$M\times S\to S$ extends to a scalar multiplication
$M_\circ\times S\to S$ whose restriction $0\times S\to S$ factors
through the zero section $0_S$ (and likewise for right modules
and bimodules).

(vi)\ \ 
Let $T$ be a topos. The notions introduced thus far for
non-pointed $T$-monoids, also admit pointed variants. Thus,
a pointed module $(S,0_S)$ is said to be {\em of finite type}
if the same holds for $S$, and $S$ is {\em finitely presented\/}
if, locally on $T$, it is the coequalizer of two morphisms
between free $\underline M$-modules of finite type.
\end{remark}

\begin{example}\label{ex_pointed-mods}
Take $T:=\Set$, and let $M$ be any monoid; then a pointed
left $M$-module is just a left $M$-module $S$ endowed with
a distinguished {\em zero element\/} $0\in S$, such that
$m\cdot 0=0$ for every $m\in M$.
A morphism $\phi:S\to S'$ of pointed left $M$-modules is just a
morphism of left $M$-modules such that $\phi(0)=0$ (and similarly
for right modules and bimodules.)

Likewise, a pointed monoid is endowed with a distinguished
{\em zero element\/}, denoted $0$ as usual, such that $0\cdot x=0$
for every $x\in M$.
\end{example}

\begin{remark}\label{rem_point-and-complete}
Let $T$ be a category as in remark \ref{rem_commutes-forgets},
and $\underline M$ a $T$-monoid.

(i)\ \  
Regardless of whether $\underline M$ is pointed or not, the category
$\underline M\Mod_{l\circ}$ is also complete and cocomplete; for
instance, if $(S,0_S)$ and $(S',0_{S'})$ are two pointed modules, the
coproduct $(S'',0_{S''}):=(S,0_S)\oplus(S',0_{S'})$ is defined by
the push-out (in the category $\underline M\Mod_l$) of the cocartesian
diagram :
$$
\xymatrix{ 0\oplus 0 \ar[rr]^-{0_S\oplus 0_{S'}} \ar[d] && S\oplus
S'\ar[d] \\ 0 \ar[rr]^{0_{S''}} && S''.}
$$
Likewise, if $\phi':S'\to S$ and $\phi'':S''\to S$ are two morphisms
in $\underline M\Mod_{l\circ}$, the fibre product $S'\times_SS''$ in
the category $\underline M\Mod_l$ is naturally pointed, and represents
the fibre product in the category of pointed modules. All this holds
also for right modules and bimodules.

(ii)\ \
The forgetful functor $\underline M\Mod_{l\circ}\to T_\circ:=1_T/T$
to the category of pointed objects of $T$, commutes with all limits,
since it is a right adjoint; it also commutes with all colimits.
This forgetful functor admits a left adjoint, that assigns to any
$\Sigma\in\Ob(T)$ the {\em free pointed $\underline M$-module\/}
$\underline M^{(\Sigma)_\circ}$. If $\underline M$ is pointed,
the latter is defined as the push-out in the cocartesian diagram
$$
\xymatrix{ 1_T\times\Sigma \ar[r] \ar[d] & M\times\Sigma \ar[d] \\
1_T \ar[r] & \underline M^{(\Sigma)_\circ}
}$$
and if $\underline M$ is not pointed, one defines it via the
equivalence of remark \ref{rem_apparent-reasons}(v) : by a
simple inspection we find that in this case
$\underline M^{(\Sigma)_\circ}=(\underline M^{(\Sigma)})_\circ$,
where $\underline M^{(\Sigma)}$ is the free (unpointed)
$\underline M$-module, as in remark \ref{rem_commutes-forgets}(iii).

Notice as well that the forgetful functors $T_\circ\to T$ and
$\underline M\Mod_{l\circ}\to\underline M\Mod_l$ both commute with
all connected colimits, hence the same also holds for the forgetful
functor $\underline M\Mod_{l\circ}\to T$. (See definition
\ref{def_filtered-cols}(vii).) The same can be repeated for right
modules and bimodules.

(iii)\ \ 
Moreover, if $\phi:S\to S'$ is any morphism in
$\underline M\Mod_{l\circ}$, we may define $\Ker\,\phi$ and
$\Coker\,\phi$ (in the category $\underline M\Mod_{l\circ}$);
namely, the kernel is the limit of the diagram
$S\xrightarrow{\phi}S'\leftarrow 0$ and the cokernel is the colimit
of the diagram $0\leftarrow S\xrightarrow{\phi}S'$. Especially, if
$S$ is a submodule of $S'$, we have a well defined quotient $S'/S$
of pointed left $\underline M$-modules. 
Furthermore, we say that a sequence of morphisms of pointed
left $\underline M$-modules :
$$
0\to S'\xrightarrow{\phi} S\xrightarrow{\psi} S''\to 0
$$
is {\em right exact}, if $\psi$ induces an isomorphism
$\Coker\,\phi\isom S''$; we say that it is {\em left exact}, if
$\phi$ induces an isomorphism $S'\isom\Ker\,\psi$, and it is {\em
short exact\/} if it is both left and right exact. (Again, all
this can be repeated also for right modules and bimodules.)
\end{remark}

\begin{example}\label{ex_rank-free-point-mod}
Take $T=\Set$, and let $M$ be a pointed or not-pointed monoid.
Then the argument from example \ref{ex_rank-mod} can be repeated
for the free pointed $M$-modules : if $\Sigma$ is any set, we have
$$
M^{(\Sigma)_\circ}\otimes_M\{1\}\isom\{1\}^{(\Sigma)_\circ}=
\Sigma_\circ
$$
where $\Sigma_\circ$ is the disjoint union of $\Sigma$ and
the final object of $\Set$ (a set with one element). Hence,
the cardinality of $\Sigma$ is an invariant, called the
{\em rank\/} of the free pointed $M$-module $M^{(\Sigma)_\circ}$,
and denoted $\rk^\circ_M M^{(\Sigma)_\circ}$.
\end{example}

\sset\subsubsection{}\label{eq_finally-flatnes}
Let $T$ be a topos, $(\underline M,0_M)$, $(\underline N,0_N)$
and $(\underline P,0_P)$ three pointed $T$-monoids, $S$, (resp.
$S'$) a pointed $(\underline M,\underline N)$-bimodule
(resp. $(\underline P,\underline N)$-bimodule); we denote
$$
\Hom_{(\underline N,0_N)_r}(S,S')
$$
the set of all morphisms of pointed right $\underline N$-modules
$S\to S'$. As usual, the presheaf
$$
\cHom_{(\underline N,0_N)_r}(S,S')
\quad :\quad
U\mapsto\Hom_{(\underline N,0_N)_{r|U}}(S_{|U},S'_{|U})
$$
(with obvious notation) is a sheaf on $(T,C_T)$, hence it
is represented by an object of $T$. Indeed, this object
is also the fibre product in the cartesian diagram :
$$
\xymatrix{
\cHom_{(\underline N,0_N)_r}(S,S') \ar[r] \ar[d] &
\cHom_{N_r}(S,S') \ar[d]^{0_S^*} \\
\cHom_{N_r}(0,0) \ar[r]^-{0_{S'*}} & \cHom_{N_r}(0,S').
}$$
Especially, $\cHom_{(\underline N,0_N)_r}(S,S')$ is naturally
a $(\underline P,\underline M)$-bimodule, and moreover,
it is pointed : its zero section represents the unique
morphism $S\to S'$ which factors through $0$.

Notice also that, for every pointed
$(\underline P,\underline M)$-bimodule $S''$, the tensor
product $S''\otimes_MS$ is naturally pointed, and as in
the non-pointed case, the functor
\set\begin{equation}\label{eq_point-tensor}
(\underline P,\underline M)\Mod_\circ\to
(\underline P,\underline N)\Mod_\circ
\quad : \quad
S''\mapsto S''\otimes_MS
\end{equation}
is left adjoint to the functor
$$
(\underline P,\underline N)\Mod_\circ\to
(\underline P,\underline M)\Mod_\circ
\quad : \quad
S'\mapsto\cHom_{(\underline N,0_N)_r}(S,S').
$$
By general nonsense, the functor \eqref{eq_point-tensor}
is right exact; especially, for any right exact sequence $T'\to
T\to T''\to 0$ of pointed $(\underline P,\underline M)$-bimodules,
the induced sequence
$$
T'\otimes_MS\to T\otimes_MS\to T''\otimes_MS\to 0
$$
is again right exact.

\begin{remark}\label{rem_point-unpointed-topos}
Suppose $\underline M$, $\underline N$ and $\underline P$ are
non-pointed $T$-monoids, $S$ is a
$(\underline M,\underline N)$-bimodule and $S''$ a
$(\underline P,\underline M)$-bimodule.

(i)\ \
If $S$ and $S''$ are pointed, one may define a tensor product
$S''\otimes_MS$ in the category
$(\underline P,\underline N)\Mod_\circ$, if one regards $S$
as a pointed
$(\underline M{}_\circ,\underline N{}_\circ)$-bimodule, and
$S''$ as a
$(\underline P{}_\circ,\underline M{}_\circ)$-bimodule as in
remark \ref{rem_apparent-reasons}(v); then one sets simply
$S''\otimes_MS:=S''\otimes_{M_\circ}S$, which is then viewed
as a pointed $(\underline P,\underline N)$-bimodule. In this
way one obtains a left adjoint to the corresponding internal
$\Hom$-functor $\cHom_N$ from pointed
$(\underline P,\underline N)$-bimodules to pointed
$(\underline P,\underline M)$-bimodules (details left to the
reader).

(ii)\ \ 
Lastly, if neither $S$ nor $S''$ is pointed, notice
the natural isomorphism :
$$
(S''\otimes_MS)_\circ\isom S''_\circ\otimes_{M_\circ}S_\circ
\qquad
\text{in the category $(\underline P,\underline N)\Mod_\circ$.}
$$
\end{remark}

\begin{definition}\label{def_pointed-flat}
In the situation of \eqref{eq_finally-flatnes}, let
$P=N:=(1_T)_\circ$, and notice that -- with these choices of
$P$ and $N$ -- a pointed $(\underline P,\underline M)$-bimodule
(resp. a pointed $(\underline M,\underline N)$-bimodule)
is just a right pointed $\underline M$-module (resp. a left
pointed $\underline M$-module), and a pointed
$(\underline P,\underline N)$-module is just a pointed object
of $T$.
\begin{enumerate}
\item
We say that $S$ is a {\em flat\/} pointed left
$\underline M$-module (or briefly, that $S$ is
{\em $\underline M$-flat}), if the functor \eqref{eq_point-tensor}
transforms short exact sequences of right pointed
$\underline M$-modules, into short exact sequences of pointed
$T$-objects. Likewise, we define flat pointed right
$\underline M$-modules.
\item
Let $\phi:\underline M\to\underline M'$ be a morphism of pointed
$T$-monoids. We say that $\phi$ is {\em flat}, if $\underline M'$
is a flat left $\underline M$-module, for the module structure
induced by $\phi$.
\end{enumerate}
\end{definition}

\begin{remark}\label{rem_flatness}
(i)\ \
In the situation of remark \ref{rem_push-pull-mond}(i),
suppose that $\underline M:=(M,0_M)$ is a pointed $T$-monoid
and $\underline N:=(N,0_N)$ a pointed $S$-monoid. By arguing
as in remark \ref{rem_push-pull-mond}(ii), we see that
$f^*\underline N:=(f^*N,f^*0_N)$ is a pointed $T$-monoid, and
$f_*\underline M:=(f_*M,f_*0_M)$ is a pointed $S$-monoid.

(ii)\ \
Likewise, if $(X,0_X)$ is a pointed left $\underline M$-module,
and $(Y,0_Y)$ a pointed left $\underline N$-module, the
$f_*(X,0_X):=(f_*X,f_*0_X)$ is a pointed $f_*\underline M$-module,
and $f^*(Y,0_Y):=(f^*Y,f^*0_Y)$ is a pointed
$f^*\underline N$-module (and likewise for right modules
and bimodules).

(iii)\ \
Also, just as in remark \ref{rem_push-pull-mond}(vii), the
associated sheaf functor $F\mapsto F^a$ transforms a presheaf
$\underline M$ of pointed monoids on $T$, into a pointed
$T$-monoid $\underline M^a$, and sends pointed left (resp.
right, resp. bi-) $\underline M$-modules to pointed left
(resp. right, resp. bi-) $\underline M^a$-modules.

(iv)\ \
Moreover, if $\phi:f^*\underline N\to\underline M$ is a morphism
of pointed $T$-monoids, then -- in view of the discussion of
\eqref{eq_finally-flatnes} -- the adjunction of remark
\ref{rem_push-pull-mond}(vi) extends to pointed modules : we
leave the details to the reader.

(v)\ \
Furthermore, in the situation of \eqref{subsec_ext-mod-byempty},
we may also define a functor
$$
j_{U!}:\underline M{}_{|U}\Mod_{l\circ}\to\underline M\Mod_{l\circ}
$$
which will be a left adjoint to $j_U^*$. Indeed, let $(X,0_X)$
be a left pointed $\underline M{}_{|U}$-module; the functor from
\eqref{subsec_ext-mod-byempty} yields a morphism $j_{U!}0_X$
of (non-pointed) $\underline M$-modules, and we define
$j_{U!}(X,0_X)$ to be the push-out (in the category
$\underline M\Mod_l$) of the diagram
$0\leftarrow j_{U!}0_X\xrightarrow{\ j_{U!}0_X }j_{U!}X$.
The latter is endowed with a natural morphism $0\to j_{U!}(X,0_X)$,
so we have a well defined pointed left $\underline M$-module.
We leave to the reader the verification that the resulting
functor, called {\em extension by zero}, is indeed left
adjoint to the restriction functor.

(vi)\ \
It is convenient to extend definition \ref{def_pointed-flat}
to non-pointed modules and monoids : namely, if $S$ is a
non-pointed left $\underline M$-module, we shall say that
$S$ is {\em flat}, if the same holds for the pointed left
$\underline M{}_\circ$-module $S_\circ$. Likewise, we say that
a morphism $\phi:\underline M\to\underline N$ of non-pointed
$T$-monoids is {\em flat}, if the same holds for $\phi_\circ$.
\end{remark}

\begin{lemma}\label{lem_platitudes}
Let $T$ be a topos, $U$ any object of $T$, and denote by
$i_*:\sC U\to T$ the inclusion functor of the complement
of $U$ in $T$ (see example {\em\ref{ex_localization-topos}(iii)}).
Let also $\underline M$, $\underline N$, $\underline P$ be
three pointed $T$-monoids. Then the following holds :
\begin{enumerate}
\item
The functor $j_{U!}$ of extension by zero is faithful,
and transforms exact sequences of pointed left
$\underline M{}_{|U}$-modules, into exact sequences of
pointed left $\underline M$-modules (and likewise for right
modules and bimodules).
\item
For every pointed $(\underline M,\underline N)$-bimodule $S$
and every pointed
$(\underline P{}_{|U},\underline M{}_{|U})$-bimodule $S'$,
the natural morphism of pointed
$(\underline P,\underline N)$-modules
$$
j_{U!}(S'\otimes_{M_{|U}}S_{|U})\to j_{U!}S'\otimes_MS
$$
is an isomorphism.
\item
If $S$ is flat pointed left $\underline M{}_{|U}$-module, then
$j_{U!}S$ is a flat pointed left $\underline M$-module (and
likewise for right modules).
\item
For every pointed $(\underline M,\underline N)$-bimodule $S$,
and every pointed $(i^*\underline P,i^*\underline M)$-bimodule
$S'$, the natural morphism of pointed
$(\underline P,\underline N)$-bimodules
$$
i_*S'\otimes_MS\to i_*(S'\otimes_{i^*M}i^*S)
$$
is an isomorphism.
\item
If $S$ is a flat pointed left $i^*\underline M$-module, then
$i_*S$ is a flat left pointed $\underline M$-module (and
likewise for right modules).
\item
If $S$ is a flat pointed left $\underline M$-module, then
$S_{|U}$ is a flat left pointed $\underline M_{|U}$-module.
\end{enumerate}
\end{lemma}
\begin{proof}(i): Let us show first that $j_{U!}$ is faithful.
Indeed, suppose that $\phi,\psi:S\to S'$ are two morphisms
of left pointed $\underline M{}_{|U}$-modules, such that
$j_{U!}\phi=j_{U!}\psi$. We need to show that $\phi=\psi$.
Let $p:S'\to S''$ be the coequalizer of $\phi$ and $\psi$;
then $j_{U!}p$ is the coequalizer of $j_{U!}\phi$ and $j_{U!}\psi$
(since $j_{U!}$ is right exact); hence we are reduced to
showing that a morphism $p:S'\to S''$ is an isomorphism
if and only if the same holds for $j_{U!}p$. This follows
from remark \ref{rem_commutes-forgets}(ii) and the following
more general :

\begin{claim}\label{cl_faithful}
Let $\phi:X\to X'$, $A\to X$, $A\to B$ be three morphisms
in $T$. Then $\phi$ is a monomorphism (resp. an epimorphism)
if and only if the same holds for the induced morphism
$\phi\amalg_AB:X\amalg_AB\to X'\amalg_AB$.
\end{claim}
\begin{pfclaim} We may assume that $T=C^\sim$ for some small
site $C:=(\cC,J)$. Then $\phi\amalg_AB=(i\phi\amalg_{iA}iB)^a$,
where $i:C^\sim\to\cC^\wedge$ is the forgetful functor.
Since the functor $F\mapsto F^a$ is exact, we are reduced
to the case where $T=\cC^\wedge$, and in this case the
assertion can be checked argumentwise, {\em i.e.} we may
assume that $T=\Set$, where the claim is obvious.
\end{pfclaim}

Next, we already know that $j_{U!}$ transforms right exact
sequences into right exact sequences. To conclude, it suffices
then to check that $j_{U!}$ transforms monomorphisms
into monomorphisms.To this aim, we apply again remark
\ref{rem_commutes-forgets}(ii) and claim \ref{cl_faithful}.

(ii) is proved by general nonsense, and (iii) is an immediate
consequence of (i) and (ii) : we leave the details to the reader.

(iv): By \eqref{subsec_up_and-down-mods}, we have
$j_U^*(i_*S'\otimes_MS)\simeq 0\otimes_{j_U^*M}j_U^*S\simeq 1_{T/U}$,
hence $i_*S'\otimes_MS\in\Ob(\sC U)$. Notice now that, for every
object $X$ of $\sC U$, the counit of adjunction $i^*i_*X\to X$
is an isomorphism (proposition \ref{prop_fullfaith-adjts}(iii));
by the triangular identities of \eqref{subsec_adj-pair}, it follows
that the same holds for the unit of adjunction $i_*X\to i_*i^*i_*X$.
Especially, the natural morphism :
$$
i_*S'\otimes_MS\to i_*i^*(i_*S'\otimes_MS)\isom
i_*(i^*i_*S'\otimes_{i^*M}i^*S)\isom i_*(S'\otimes_{i^*M}i^*S).
$$
is an isomorphism. The latter is the morphism of assertion (iv).

(v) follows easily from (iv) and its proof.

(vi): In view of (i), it suffices to show that the functor
$S'\mapsto j_{U!}(S'\otimes_{M_{|U}}S_{|U})$ transforms exact
sequences into exact sequences. The latter follows easily from
(ii).
\end{proof}

\begin{proposition}\label{prop_finally-flat}
Let $\bP(T,\underline M,S)$ be the property : ``$S$ is a
flat pointed left $\underline M$-module'' (for a monoid
$\underline M$ on a topos $T$). Then $\bP$ can be checked
on stalks. (See definition {\em\ref{def_T-point}(v)}.)
\end{proposition}
\begin{proof} Suppose first that $S_\xi$ is a flat left
$\underline M{}_\xi$-module for every $\xi$ in a
conservative set of $T$-points; let $\phi:X\to X'$
be a monomorphism of pointed right $\underline M$-modules;
by \eqref{subsec_up_and-down-mods} we have a natural isomorphism
$$
(\phi\otimes_MS)_\xi\isom\phi_\xi\otimes_{M_\xi}S_\xi
$$
in the category of pointed sets, and our assumption implies
that these morphisms are monomorphisms. Since an arbitrary
product of monomorphisms is a monomorphism, remark
\ref{rem_cofinal}(ii) shows that $\phi\otimes_MS$ is
also a monomorphism, whence the contention.

Next, suppose that $S$ is a flat pointed left
$\underline M$-module. We have to show that the functor
\set\begin{equation}\label{eq_yaik}
S'\mapsto S'\otimes_{M_\xi}S_\xi
\end{equation}
from pointed right $\underline M_\xi$-modules to pointed
sets, preserves monomorphisms.

However, let $(U,\xi_U,\omega_U)$ be any lifting of $\xi$
(see \eqref{subsec_blift}); in view of
\eqref{subsec_up_and-down-mods}, we have 
$$
P_U(S'):=(\xi_U^*\xi_{U*}S')\otimes_{A_\xi}S_\xi\simeq
(\xi_U^*\xi_{U*}S')\otimes_{\xi^*_UA_{|U}}\xi^*_US_{|U}
\simeq\xi_U^*(\xi_{U*}S'\otimes_{A_{|U}}S_{|U})
$$
and then lemma \ref{lem_platitudes}(vi) implies that the
functor $S'\to P_U(S')$ preserves monomorphisms. By lemma
\ref{lem_super-cazzuto}, remark \ref{rem_Nbd-as-Fib}(ii) and
\eqref{eq_same-neighbor}, the functor \eqref{eq_yaik}
is a filtered colimit of such functors $P_U$, hence it
preserves monomorphisms as well.
\end{proof}

\sset\subsubsection{}\label{subsec_wish-now}
We wish now to introduce a few notions that pertain to the
special class of commutative $T$-monoids. When $T=\Set$,
these notions are well known, and we wish to explain quickly
that they generalize without problems, to arbitrary topoi.

To begin with, for every category $T$ as in example
\ref{ex_stupid-prior}(i), we denote by $\Mnd_T$ (resp.
$\Mnd_{T\circ}$) the category of commutative non-pointed
(resp. pointed) $T$-monoids; in case $T=\Set$, we shall usually
drop the subcript, and write just $\Mnd$ (resp. $\Mnd_\circ$).
Notice that if $\underline M$ is any (pointed or not pointed)
commutative $T$-monoid, every left or right $\underline M$-module
is a $\underline M$-bimodule in a natural way, hence we shall
denote indifferently by $\underline M\Mod$ (resp.
$\underline M\Mod{}_\circ$) the category of non-pointed (resp.
pointed) left or right $\underline M$-modules. 

The following lemma is a special case of a result that holds more
generally, for every "algebraic theory" in the sense of
\cite[Def.3.3.1]{BorII} (see \cite[Prop.3.4.1, Prop.3.4.2]{BorII}).

\begin{lemma}\label{lem_forget-me-not}
Let $T$ be a topos. We have :
\begin{enumerate}
\item
The category $\Mnd_T$ admits arbitrary limits and colimits.
\item
In the category $\Mnd_T$, filtered colimits commute with all finite
limits.
\item
The forgetful functor $\iota:\Mnd_T\to T$ that assigns to a monoid
its underlying object of $T$, commutes with all limits, and with all
filtered colimits.
\end{enumerate}
\end{lemma}
\begin{proof} (iii): Commutation with limits holds because $\iota$
admits a left adjoint : namely, to an object $\Sigma$ of $T$ one
assigns the {\em free monoid\/} $\N^{(\Sigma)}_T$ generated by
$\Sigma$, defined as the sheaf associated with the presheaf of
monoids
$$
U\mapsto\N^{(\Sigma(U))}
\qquad
\text{for every $U\in\Ob(T)$}
$$
where $\N$ is the additive monoid of natural numbers
(see remark \ref{rem_push-pull-mond}(vii)). One verifies
easily that this $T$-monoid represents the functor
$$
\underline M\mapsto\Hom_T(\Sigma,\underline M)
\qquad
\Mnd_T\to\Set.
$$
Moreover, if $I$ is any small category, and $F:I\to\Mnd_T$
any functor, one checks easily that the limit of $\iota\circ F$
can be endowed with a unique composition law (indeed, the limit
of the composition laws of the monoids $F_i$), such that
the resulting monoid represents the limit of $F$.

A similar argument also shows that $\Mnd_T$ admits arbitrary filtered
colimits, and that $\iota$ commutes with filtered colimits. It is
likewise easy to show that the product of two $T$-monoids
$\underline M$ and $\underline N$ is also the coproduct of
$\underline M$ and $\underline N$. To complete the proof of (i), it
suffices therefore to show that any two maps
$f,g:\underline M\to\underline N$ admit a coequalizer; the latter
is obtained as the coequalizer $\underline N'$ (in the category $T$)
of the two morphisms :
$$
\xymatrix{ \underline M\times\underline N 
\ar@<-.5ex>[rr]_-{\mu_N\circ(g\times\one_N)}
\ar@<.5ex>[rr]^-{\mu_N\circ(f\times\one_N)} & & \underline N.
}$$
We leave to the reader the verification that the composition law
of $N$ descends to a (necessarily unique) composition law on
$\underline N'$.

(ii) follows from (iii) and the fact that the same assertion
holds in $T$ (remark \ref{rem_summarized}(i)).
\end{proof}

\begin{example}\label{ex_forrinst}
(i)\ \
For instance, if $T=\Set$, the product $M_1\times M_2$ of any two
commutative monoids is representable in $\Mnd$; its underlying set
is the cartesian product of $M_1$ and $M_2$, and the composition
law is the obvious one.

(ii)\ \
As usual, the kernel $\Ker\,\phi$ (resp. cokernel $\Coker\,\phi$)
of a map of $T$-monoids $\phi:\underline M\to\underline N$ is
defined as the fibre product (resp. push-out) of the diagram of
$T$-monoids
$$
\underline M\xrightarrow{\ \phi\ }
\underline N\leftarrow\underline 1{}_T
\qquad
\text{(resp. $\underline 1{}_T\leftarrow\underline M
\xrightarrow{\ \phi\ }\underline N$)}.
$$
Especially, if $\underline M\subset\underline N$, one defines
in this way the quotient $\underline N/\underline M$.

(iii)\ \
Also, if $T=\Set$, and $\phi_1:M\to M_1$, $\phi_2:M\to M_2$ are
two maps in $\Mnd$, the push-out $M_1\amalg_MM_2$ can be described
as follows. As a set, it is the quotient $(M_1\times M_2)/\!\!\sim$,
where $\sim$ denotes the minimal equivalence relation such that
$$
(m_1,m_2\cdot\phi_2(m))\sim(m_1\cdot\phi_1(m),m_2) \qquad\text{for
every $m\in M$, $m_1\in M_1$, $m_2\in M_2$}
$$
and the composition law is the unique one such that the projection
$M_1\times M_2\to M_1\amalg_MM_2$ is a map of monoids. We deduce
the following :
\end{example}

\begin{lemma}\label{lem_special-p-out}
Let $G$ be an abelian group. The following holds :
\begin{enumerate}
\item
If $\phi:M\to N$ and $\psi:M\to G$ are two morphisms of monoids
(in the topos $T=\Set$), $G\amalg_MN$ is the quotient
$(G\times N)/\!\!\approx$, where $\approx$ is the equivalence
relation such that :
$$
(g,n)\approx(g',n')
\quad\Leftrightarrow\quad
\text{$(\psi(a)\cdot g,\phi(b)\cdot n)=
(\psi(b)\cdot g',\phi(a)\cdot n')$ for some $a,b\in M$}.
$$
\item
If $\phi:G\to M$ and $\psi:G\to N$ are two morphisms of monoids,
the set underlying $M\amalg_GN$ is the set-theoretic quotient
$(M\times N)/G$ for the $G$-action defined via $(\phi,\psi^{-1})$.
\item
Especially, if $M$ is a monoid and $G$ is a submonoid of $M$, then
the set underlying $M/G$ is the set-theoretic quotient of $M$ by the
translation action of\/ $G$.
\end{enumerate}
\end{lemma}
\begin{proof}(i): One checks easily that the relation $\approx$
thus defined is transitive. Let $\sim$ be the equivlence relation
defined as in example \ref{ex_forrinst}(iii). Clearly :
$$
(g,n\cdot\psi(m))\approx(g\cdot\phi(m),n) \qquad\text{for every
$g\in G$, $n\in N$ and $m\in M$}
$$
hence $(g,n)\sim(g',n')$ implies $(g,n)\approx(g',n')$. Conversely,
suppose that $(\psi(a)\cdot g,\phi(b)\cdot n)= (\psi(b)\cdot
g',\phi(a)\cdot n')$ for some $g\in G$, $n\in N$ and $a,b\in M$.
Then :
$$
(g,n)=(g,\phi(a)\cdot\phi(a)^{-1}\cdot n)\sim (\psi(a)\cdot
g,\phi(a)^{-1}\cdot n)= (\psi(b)\cdot g',\phi(a)^{-1}\cdot n)
$$
as well as : $(g',n')=(g',\phi(b)\cdot\phi(a)^{-1}\cdot n)\sim
(\psi(b)\cdot g',\phi(a)^{-1}\cdot n)$. Hence $(g,n)\sim(g',n')$
and the claim follows.

(ii) follows directly from example \ref{ex_forrinst}(iii), and (iii)
is a special case of (ii).
\end{proof}

\sset\subsubsection{}\label{subsec_regular-mon}
Let $T$ be a topos.
For any $T$-ring $\underline R$, we let $\underline R\Mod$ be
the category of $\underline R$-modules (defined in the usual way);
especially, we may consider the $T$-ring $\Z_T$ (the constant
sheaf with value $\Z$ : see \eqref{subsec_glob-sections}).
Then $\Z_T\Mod$ is the category of abelian $T$-groups.
The forgetful functor $\Z_T\Mod\to\Mnd_T$ admits a right adjoint :
$$
\Mnd_T\to\Z_T\Mod
\quad:\quad
\underline M\mapsto\underline M^\times.
$$
The latter can be defined as the fibre product in the cartesian
diagram :
$$
\xymatrix{ \underline M^\times \ar[r] \ar[d] & 
\underline M\times\underline M \ar[d]^{\mu_M} \\
1_T \ar[r]^-{1_M} & M.
}$$
For $i=1,2$, let $p_i:\underline M\times\underline M\to\underline M$
be the projections, and $p'_i:\underline M^\times\to\underline M$
the restriction of $p_i$; for every $U\in\Ob(T)$, the image of
$p'_i(U):\underline M^\times(U)\to\underline M(U)$ consists of all
sections $x$ which are {\em invertible}, {\em i.e.} for which there
exists $y\in\underline M(U)$ such that $\mu_M(x,y)=1_M$. It is
easily seen that such inverse is unique, hence $p'_i$ is a
monomorphism, $p'_1$ and $p'_2$ define the same subobject of
$\underline M$, and this subobject $\underline M^\times$ is the
largest abelian $T$-group contained in $\underline M$. We say
that $\underline M$ is {\em sharp}, if $\underline M^\times=1_T$.
The inclusion functor, from the full subcategory of sharp $T$-monoids,
to $\Mnd_T$, admits a left adjoint
$$
\underline M\mapsto\underline M^\sharp:=
\underline M/\underline M^\times.
$$
We call $\underline M^\sharp$ the {\em sharpening\/} of
$\underline M$.

\sset\subsubsection{}\label{subsec_localize-mons}
Let $\underline S$ be a submonoid of a commutative $T$-monoid
$\underline M$, and $F_S:\Mnd_T\to\Set$ the functor that assigns
to any commutative $T$-monoid $\underline N$ the set of all
morphisms $f:\underline M\to\underline N$ such that
$f(\underline S)\subset\underline N^\times$. We claim that $F_S$
is representable by a $T$-monoid $\underline S^{-1}\underline M$.

In case $T=\Set$, one may realize $\underline S^{-1}\underline M$
as the quotient $(\underline S\times\underline M)/\!\!\sim$ for
the equivalence relation such that $(s_1,x_1)\sim(s_2,x_2)$ if
and only if there exists $t\in\underline S$ such that
$ts_1x_2=ts_2x_1$. The composition law of
$\underline S^{-1}\underline M$ is the obvious one; then the class
of a pair $(s,x)$ is denoted naturally by $s^{-1}x$. This construction
can be repeated on a general topos : letting
$X:=\underline S\times\underline M$, the foregoing equivalence
relation can be encoded as the equalizer $R$ of two maps
$X\times X\times S\to\underline M$, and the
quotient under this equivalence relation shall be represented
by the coequalizer of two other maps $R\to X$; the reader may
spell out the details, if he wishes. Equivalently,
$\underline S^{-1}\underline M$ can be realized as the sheaf
on $(T,C_T)$ associated with the presheaf :
$$
T\to\Mnd
\quad :\quad
U\mapsto\underline S(U)^{-1}\underline M(U)
$$
(see remark \ref{rem_push-pull-mond}(vii)). The natural morphism
$\underline M\to\underline S^{-1}\underline M$ is called the
{\em localization map}. For $T=\Set$, and $f\in M$ any element,
we shall also use the standard notation :
$$
M_f:=S^{-1}_fM \qquad \text{where $S_f:=\{f^n~|~n\in\N\}$}.
$$

\begin{lemma}\label{lem_localize}
Let $f_1:\underline M\to\underline N{}_1$ and
$f_2:\underline M\to\underline N{}_2$ be morphisms of $T$-monoids,
$\underline S\subset\underline M$,
$\underline S{}_i\subset\underline N{}_i$ ($i=1,2$) three submonoids,
such that $f_i(\underline S)\subset\underline S{}_i$ for $i=1,2$.
Then the natural morphism :
$$
(\underline S{}_1\cdot\underline S{}_2)^{-1}
(\underline N{}_1\amalg_{\underline M}\underline N{}_2)\to
\underline S{}_1^{-1}\underline N{}_1
\amalg_{\underline S^{-1}\underline M}
\underline S{}_2^{-1}\underline N{}_2
$$
is an isomorphism.
\end{lemma}
\begin{proof} One checks easily that both these $T$-monoids
represent the functor $\Mnd_T\to\Set$ that assigns to any
$T$-monoid $\underline P$ the pairs of morphisms $(g_1,g_2)$
where $g_i:\underline N{}_i\to\underline P$ satisfies
$g_i(\underline S{}_i)\subset\underline P^\times$, for
$i=1,2$, and $g_1\circ f_1=g_2\circ f_2$. The details are
left to the reader.
\end{proof}

\sset\subsubsection{}\label{subsec_int-mds}
The forgetful functor $\Z\Mod_T\to\Mnd_T$ from abelian $T$-groups
to commutative $T$-monoids, admits a left adjoint
$$
\underline M\mapsto\underline M^\gp:=\underline M^{-1}\underline M.
$$
A commutative $T$-monoid $\underline M$ is said to be {\em integral\/}
if the unit of adjunction $\underline M\to\underline M^\gp$ is a
monomorphism. The functor $\underline M\mapsto\underline M^\gp$
commutes with all colimits, since all left adjoints do; it does
not commute with arbitrary limits (see example \ref{ex_regular}(v)).

We denote by $\IntMnd_T$ the full subcategory of $\Mnd_T$
consisting of all integral monoids; when $T=\Set$, we omit
the subscript, and write just $\IntMnd$. The natural inclusion
$\iota:\IntMnd_T\to\Mnd_T$ admits a left adjoint :
$$
\Mnd_T\to\IntMnd_T
\quad :\quad
\underline M\mapsto\underline M^\intg.
$$
Namely, $\underline M^\intg$ is the image (in the category
$T$) of the unit of adjunction
$\underline M\to\underline M^\gp$. It follows easily that
the category $\IntMnd_T$ is cocomplete, since the colimit
of a family $(\underline M{}_\lambda~|~\lambda\in\Lambda)$
of integral monoids is represented by
$$
(\colim_{\lambda\in\Lambda}\iota\underline M{}_\lambda)^\intg.
$$
Likewise, $\IntMnd_T$ is complete, and limits commute with the
forgetful functor to $T$; to check this, it suffices to show that
$$
L:=\lim_{\lambda\in\Lambda}\iota(\underline M_\lambda)
$$
is integral. However, by lemma \ref{lem_forget-me-not}(iii) we
have $L\subset\prod_{\lambda\in\Lambda}\underline M_\lambda\subset
\prod_{\lambda\in\Lambda}\underline M_\lambda^\gp$, whence the
claim.

\begin{example}\label{ex_regular}
(i)\ \ 
Take $T=\Set$; if $M$ is any monoid, and $a\in M$ is any element,
we say that $a$ is {\em regular}, if the map $M\to M$ given by
the rule : $x\mapsto a\cdot x$ is injective. It is easily seen
that $M$ is integral if and only if every element of $M$ is regular.

(ii)\ \
For an arbitrary topos $T$, notice that the $T$-monoid
$\underline G^a$ associated with a presheaf of groups
$\underline G$ on $T$, is a $T$-group : indeed, the
condition $\underline G^\times=\underline G$ implies
$\underline (G^a)^\times=\underline G^a$, since the
functor $F\mapsto F^a$ is exact. More precisely, for
every presheaf $\underline M$ of monoids on $T$, we have
a natural isomorphism :
$$
(\underline M^\gp)^a\isom(\underline M^a)^\gp
\qquad
\text{for every $T$-monoid $\underline M$}
$$
since both functors are left adjoint to the forgetful
functor from $T$-groups to presheaves of monoids on $T$.

(iii)\ \
It follows from (ii) that a $T$-monoid $\underline M$ is integral
if and only if $\underline M(U)$ is an integral monoid, for every
$U\in\Ob(T)$. Indeed, if $\underline M$ is integral, then
$\underline M(U)\subset\underline M^\gp(U)$ for every such
$U$, so $\underline M(U)$ is integral. Conversely, by definition
$\underline M^\gp$ is the sheaf associated with the presheaf
$U\mapsto\underline M(U)^\gp$; now, if $\underline M(U)$
is integral, we have $\underline M(U)\subset\underline M(U)^\gp$,
and consequently $\underline M\subset\underline M^\gp$, since
the functor $F\mapsto F^a$ is exact.

(iv)\ \
We also deduce from (ii) that the functor
$\underline M\mapsto\underline M^a$ sends presheaves of
integral monoids, to integral $T$-monoids. Therefore we
have a natural isomorphism :
\set\begin{equation}\label{eq_already-intg}
(\underline M^\intg)^a\isom(\underline M^a)^\intg
\end{equation}
as both functors are left adjoint to the forgetful
functor from integral $T$-monoids, to presheaves of monoids
on $T$. In the same vein, it is easily seen that the forgetful
functor $\IntMnd_T\to T$ commutes with filtered colimits : indeed,
\eqref{eq_already-intg} and lemma \ref{lem_forget-me-not}(iii)
reduce the assertion to showing that the colimit of a filtered
system of presheaves of integral monoids is integral, which can
be verified directly.

(v)\ \
Take $T=\Set$, and let $\phi:M\to N$ be an injective map of monoids;
if $N$ (hence $M$) is integral, one sees easily that the induced map
$\phi^\gp:M^\gp\to N^\gp$ is also injective. This may fail, when
$N$ is not integral : for instance, if $M$ is any integral monoid,
and $N:=M_\circ$ is the pointed monoid associated with $M$ as in
remark \ref{rem_apparent-reasons}(iv), then for the natural inclusion
$i:M\to M_\circ$ we have $i^\gp=0$, since $(M_\circ)^\gp=\{1\}$.
\end{example}

\begin{lemma}\label{lem_integral-quot}
Let $T$ be a topos with enough points, $\underline M$ an integral
$T$-monoid, and $\underline N\subset\underline M$ a $T$-submonoid.
Then $\underline M/\underline N$ is an integral $T$-monoid.
\end{lemma}
\begin{proof} In light of example \ref{ex_regular}(iv), we are
reduced to the case where $T=\Set$. Moreover, since the natural
morphism
$\underline M/\underline N\to\underline N^{-1}M/\underline N^\gp$
is an isomorphism, we may assume that $\underline N$ is an
abelian group. Now, notice that
$(\underline M/\underline N)^\gp=\underline M^\gp/\underline N$
since the functor $P\mapsto P^\gp$ commutes with colimits. On
the other hand, $\underline M/\underline N$ is the set-theoretic
quotient of
$\underline M$ by the translation action of $\underline N$
(lemma \ref{lem_special-p-out}(iii)). This shows that the unit
of adjunction
$\underline M/\underline N\to(\underline M/\underline N)^\gp$
is injective, as required.
\end{proof}

\sset\subsubsection{}
Let $M$ be an integral monoid. Classically, one says that $M$
is {\em saturated}, if we have :
$$
M=\{a\in M^\gp~|~a^n\in M \text{ for some integer $n>0$}\}.
$$
In order to globalize the class of saturated monoid to
arbitrary topoi, we make the following :

\begin{definition}\label{def_exact-phi}
Let $T$ be a topos, $\phi:\underline M\to\underline N$ a
morphism of integral $T$-monoids.
\begin{enumerate}
\item
We say that $\phi$ is {\em exact\/} if the diagram of
commutative $T$-monoids
$$
\cD_\phi \quad : \quad
{\diagram
\underline M \ar[r]^-\phi \ar[d] & \underline N \ar[d] \\
\underline M^\gp \ar[r]^-{\phi^\gp} & \underline N^\gp 
\enddiagram}$$
is cartesian (where the vertical arrows are the natural
morphisms).
\item
For any integer $k>0$, the {\em $k$-Frobenius map\/} of
$\underline M$ is the endomorphism $\bek_M$ of $\underline M$ 
given by the rule : $x\mapsto x^k$ for every $U\in\Ob(T)$ and
every $x\in\underline M(U)$. We say that $\underline M$ is
{\em $k$-saturated}, if $\bek_M$ is an exact morphism.
\item
We say that $\underline M$ is {\em saturated}, if $\underline M$
is integral and $k$-saturated for every integer $k>0$.
\end{enumerate}
\end{definition}

We denote by $\SatMnd_T$ the full subcategory of $\IntMnd_T$
whose objects are the saturated $T$-monoids. As usual, when
$T=\Set$, we shall drop the subscript, and just write
$\SatMnd$ for this category.
The above definition (and several of the related results in section
\ref{subsec_intg-monnd}) is borrowed from \cite{Tsu}. 

\begin{remark}\label{rem_satura}
(i)\ \
Clearly, when $T=\Set$, definition \ref{def_exact-phi}(iii)
recovers the classical notion of saturated monoid. Again, for
usual monoids, it is easily seen that the forgetful functor
$\SatMnd\to\IntMnd$ admits a left adjoint, that assigns to
any integral monoid $M$ its {\em saturation\/} $M^\sat$.
The latter is the monoid consisting of all elements $x\in M^\gp$
such that $x^k\in M$ for some integer $k>0$; especially, the
torsion subgroup of $M^\gp$ is always contained in $M^\sat$.
The easy verification is left to the reader. Clearly, $M$
is saturated if and only if $M=M^\sat$. More generally, the
unit of adjunction $M\to M^\sat$ is just the inclusion map.

(ii)\ \
For a general topos $T$, and a morphism $\phi$ as in definition
\ref{def_exact-phi}(i), notice that $\phi$ is exact if and only
if the induced map of monoids
$\phi(U):\underline M(U)\to\underline N(U)$ is exact for every
$U\in\Ob(T)$. Indeed, if $\cD_\phi$ is cartesian, then the same
holds for the induced diagram $\cD_\phi(U)$ of monoids; since
the natural map $\underline M(U)^\gp\to\underline M^\gp(U)$ is
injective (and likewise for $\underline N$), it follows easily
that the diagram of monoids $\cD_{\phi(U)}$ is cartesian, {\em
i.e.} $\phi(U)$ is exact.
For the converse, notice that $\cD_\phi$ is of the form
$(h\cD_\phi)^a$, where $h:T\to T^\wedge$ is the Yoneda
embedding, and $F\mapsto F^a$ denotes the associated sheaf
functor $T^\wedge\to(T,C_T)^\sim=T$; the assumption means that
$h\cD$ is a cartesian diagram in $T^\wedge$, hence $\cD$ is
exact in $T$, since the associated sheaf functor is exact.

(iii)\ \
Example \ref{ex_regular}(iii) and (ii) imply that a $T$-monoid
$\underline M$ is saturated, if and only if $\underline M(U)$
is a saturated monoid, for every $U\in\Ob(T)$.
We also remark that, in view of example \ref{ex_regular}(ii),
the functor $F\mapsto F^a$ takes presheaves of $k$-saturated
(resp. saturated) monoids, to $k$-saturated (resp. saturated)
$T$-monoids : indeed, if $\eta:\underline M\to\underline M^\gp$
is the unit of adjunction for a presheaf of monoids $\underline M$,
then
$\eta^a:\underline M^a\to(\underline M^\gp)^a=(\underline M^a)^\gp$
is the unit of adjunction for the associated $T$-monoid, hence
it is clear the functor $F\mapsto F^a$ preserves exact
morphisms.

(iv)\ \
It follows easily that the inclusion functor
$\SatMnd_T\to\IntMnd_T$ admits a left adjoint, namely the
functor
$$
\IntMnd_T\to\SatMnd_T
\quad :\quad
\underline M\mapsto\underline M^\sat
$$
that assigns to $\underline M$ the sheaf associated with the
presheaf $U\to\underline M'(U):=\underline M(U)^\sat$ on $T$
(notice that the functor $\underline M\mapsto\underline M'$
from presheaves of integral monoids, to presheaves of saturated
monoids, is left adjoint to the inclusion functor).
Just as in example \ref{ex_regular}(iv), we deduce a natural
isomorphism
\set\begin{equation}\label{eq_saturat}
(\underline M^\sat)^a\isom(\underline M^a)^\sat
\qquad
\text{for every $T$-monoid $\underline M$}
\end{equation}
since both functors are left adjoint to the forgetful functor
from $\SatMnd_T$, to presheaves of integral monoids on $T$.

(v)\ \ 
By the usual general nonsense, the saturation functor commutes
with all colimits.
Moreover, the considerations of \eqref{subsec_int-mds} can be
repeated for saturated monoids : first, the category $\SatMnd_T$ is
cocomplete, and arguing as in example \ref{ex_regular}(iv),
one checks that filtered colimits commute with the forgetful functor
$\SatMnd_T\to T$; next, if $F:\Lambda\to\SatMnd_T$ is a functor
from a small category $\Lambda$, then for each integer $k>0$,
the induced diagram of integral monoids
$$
\displaystyle\lim_\Lambda\cD_{\bek_F}
\quad :\quad
{\diagram
\displaystyle\lim_\Lambda F \ar[rr]^-{\lim_\Lambda\bek_F} \ar[d] & &
\displaystyle\lim_\Lambda F \ar[d] \\
\displaystyle\lim_\Lambda F^\gp \ar[rr]^-{\lim_\Lambda\bek_F^\gp} & & 
\displaystyle\lim_\Lambda F^\gp
\enddiagram}$$
is cartesian; since the natural morphism
$$
(\lim_\Lambda F)^\gp\to\lim_\Lambda F^\gp
$$
is a monomorphism, it follows easily that the limit of $F$
is saturated, hence $\SatMnd_T$ is complete, and furthermore
all limits commute with the forgetful functor to $T$.
\end{remark}

\sset\subsubsection{}\label{subsec_pts-of-top}
In view of remark \ref{rem_push-pull-mond}(i,ii,iii), a
morphism of topoi $f:T\to S$ induces functors :
\set\begin{equation}\label{eq_ind-ooon-mnds}
f_*:\Mnd_T\to\Mnd_S
\qquad
f^*:\Mnd_S\to\Mnd_T
\end{equation}
and one verifies easily that \eqref{eq_ind-ooon-mnds} is an adjoint
pair of functors.

\begin{lemma}\label{lem_T-satura}
Let $f:T\to S$ be a morphism of topoi, $\underline M$ an
$S$-monoid. We have :
\begin{enumerate}
\item
If $\underline M$ is integral (resp. saturated), $f^*\underline M$
is an integral (resp. saturated) $T$-monoid.
\item
More precisely, there is a natural isomorphism :
$$
f^*(\underline M^\intg)\isom(f^*\underline M)^\intg
\qquad
\text{(resp. $f^*(\underline M^\sat)\isom(f^*\underline M)^\sat$,
\ \  if $\underline M$ is integral)}.
$$
\item
If $\phi$ is an exact morphism of integral $S$-monoids, then
$f^*\phi$ is an exact morphism of integral $T$-monoids.
\end{enumerate}
\end{lemma}
\begin{proof} To begin with, notice that the adjoint
pair $(f^*,f_*)$ of \eqref{eq_ind-ooon-mnds} restricts
to a corresponding adjoint pair of functors between the
categories of abelian $T$-groups and abelian $S$-groups
(since the condition $\underline G=\underline G^\times$
for monoids, is preserved by any left exact functor).

There follows a natural isomorphism :
$$
(f^*\underline M)^\gp\isom f^*(\underline M^\gp)
\qquad
\text{for every $S$-monoid $\underline M$}
$$
since both functors are left adjoint to the functor
$f_*$ from abelian $T$-groups to $S$-monoids. Now, if
$\underline M$ is an integral $S$-module, and
$\eta:\underline M\to\underline M^\gp$ is the unit of
adjunction, it is easily seen that
$f^*\eta:f^*\underline M\to(f^*\underline M)^\gp$ is
also the unit of adjunction.
From this and proposition \ref{prop_half-dominates}(iii),
we deduce the assertion concerning $f^*(\underline M^\intg)$.

By the same token, we get assertion (iii) of the lemma.
Especially, if $\underline M$ is saturated, then the
same holds for $f^*\underline M$. The assertion concerning
$f^*(\underline M^\sat)$ follows by the usual argument.
\end{proof}

\begin{lemma}\label{lem_concern}
{\em(i)}\ \
The functor $f^*$ of \eqref{eq_ind-ooon-mnds}
commutes with all finite limits and all colimits.

{\em(ii)}\ \
Let $\bP(T,M)$ be the property ``$M$ is an integral (resp.
saturated) $T$-monoid''. Then $\bP$ can be checked on every
conservative set of morphisms of topoi (see definition
{\em\ref{def_T-point}(v)}.)
\end{lemma}
\begin{proof}(i): Concerning finite limits, in light of
lemma \ref{lem_forget-me-not}(iii) we are reduced to the
assertion that $f^*:S\to T$ is left exact, which holds by
definition. Next $f^*$ commutes with colimits, because it
is a left adjoint.

(ii): A $T$-monoid $\underline M$ is integral if and only if the
unit of adjunction $\eta:\underline M\to\underline M^\intg$
is an isomorphism. However,
$f^*(\underline M^\intg)\isom(f^*\underline M)^\intg$, by lemma
\ref{lem_T-satura}(ii), and
$f^*\eta:\underline M\to(f^*\underline M)^\intg$
is the unit of adjunction. The assertion is an immediate
consequence. The same argument applies as well to saturated
$T$-monoids.
\end{proof}

\begin{example}\label{ex_constant-mond}
(i)\ \
For instance, the unique morphism of topoi $\Gamma:T\to\Set$
(see proposition \ref{prop_glob-sections}(i,ii)) induces a pair
of adjoint functors :
\set\begin{equation}\label{eq_const-mon-top}
\Mnd_T\to\Mnd\ :\
M\mapsto\Gamma(T,M)\qquad\text{and}\qquad \Mnd\to\Mnd_T\ :\
P\mapsto P_T
\end{equation}
where $P_T$ is the constant sheaf of monoids on $(T,C_T)$ with value
$P$.

(ii)\ \
Specializing lemma \ref{lem_T-satura}(ii) to this adjoint pair,
we obtain natural isomorphisms :
\set\begin{equation}\label{eq_nat-int-glob}
(M_T)^\intg\isom(M^\intg)_T \qquad (M_T)^\sat\isom(M^\sat)_T
\end{equation}
of functors $\Mnd\to\IntMnd_T$ and $\IntMnd\to\SatMnd_T$.
Especially, if $M$ is an integral (resp. saturated) monoid,
then the constant $T$-monoid $M_T$ is integral (resp. saturated).

(iii)\ \
If $\xi$ is any $T$-point, notice also that the stalk $M_{T,\xi}$
is isomorphic to $M$, since $\xi$ is a section of $\Gamma:T\to\Set$.
\end{example}

\sset\subsubsection{}\label{subsec_mon-to-algs}
Let $T$ be a topos, $\underline R$ a $T$-ring. We have a forgetful
functor $\underline R\Alg\to\Mnd_T$ that assigns to a (unital, commutative)
$\underline R$-algebra $(\underline A,+,\cdot,1_A)$ its multiplicative
$T$-monoid $(\underline A,\cdot)$. If $T=\Set$, this functor admits
a left adjoint $\Mnd\to R\Alg$ : $M\mapsto R[M]$. Explicitly,
$R[M]=\bigoplus_{x\in M}xR$, and the multiplication law is uniquely
determined by the rule :
$$
xa\cdot yb:=(x\cdot y)ab
\qquad
\text{for every $x,y\in M$ and $a,b\in R$}.
$$
For a general topos $T$, the above construction globalizes to
give a left adjoint
\set\begin{equation}\label{eq_another-one}
\Mnd_T\to\underline R\Alg
\quad : \quad
\underline M\mapsto\underline R[\underline M].
\end{equation}
The latter is the sheaf on $(T,C_T)$ associated with the presheaf
$U\mapsto\underline R(U)[\underline M(U)]$, for every $U\in\Ob(T)$.
The functor \eqref{eq_another-one} commutes with arbitrary colimits
(since it is a left adjoint); especially, if
$\underline M\to\underline M{}_1$ and
$\underline M\to\underline M{}_2$ are two morphisms of monoids,
we have a natural identification :
\set\begin{equation}\label{eq_push-out-tensor}
\underline R[\underline M{}_1\amalg_{\underline M}\underline M{}_2]\isom
\underline R[\underline M_1]\otimes_{\underline R[\underline M]}
\underline R[\underline M_2].
\end{equation}
By inspecting the universal properties, we also get a natural
isomorphism :
\set\begin{equation}\label{eq_loc-monds}
\underline S^{-1}\underline R[\underline M]\isom
\underline R[\underline S^{-1}\underline M]
\end{equation}
for every monoid $\underline M$ and every submonoid
$\underline S\subset\underline M$.

\sset\subsubsection{}
Likewise, if $\underline M$ is any $T$-monoid, let
$\underline R[\underline M]\Mod$ denote as usual the category
of modules over the $T$-ring $\underline R[\underline M]$;
we have a forgetful functor
$\underline R[\underline M]\Mod\to\underline M\Mod$. When
$T=\Set$, this functor admits a left adjoint
$M\Mod\to R[M]\Mod$ : $S\mapsto R[S]$. Explicitly, $R[S]$ is
the free $R$-module with basis given by $S$, and the
$R[M]$-module structure on $R[S]$ is determined by the rule:
$$
xa\cdot sb:=\mu_S(x,s)ab
\qquad
\text{for every $x\in M$, $s\in S$ and $a,b\in R$}.
$$
For a general topos $T$, this construction globalizes to
give a left adjoint
$$
\underline M\Mod\to\underline R[\underline M]\Mod
\quad : \quad
(S,\mu_S)\mapsto\underline R[S]
$$
which is defined as the sheaf associated with the presheaf
$U\mapsto\underline R(U)[S(U)]$ in $T^\wedge$.

\subsection{Torsors on a topos}
In this section we discuss torsors over (not necessarily
abelian) group objects of a topos. Then we explain some
basic notions concerning the points of the \'etale and
Zariski topoi of a scheme, and we conclude with the proof
of Hilbert's theorem 90 (lemma \ref{lem_Hilbert90}(iv)).

\begin{definition}\label{def_A-torsors}
Let $T$ be a topos, and $G$ a $T$-group.
\begin{enumerate}
\item
A {\em left $G$-torsor\/} is a left $G$-module $(X,\mu_X)$,
inducing an isomorphism
$$
(\mu_X,p_X):G\times X\to X\times X
$$
(where $p_X:G\times X\to X$ is the natural projection) and such
that there exists a covering morphism $U\to 1_T$ in $T$ for which
$X(U)\neq\emptyset$. This is the same as saying that the unique
morphism $X\to 1_T$ is an epimorphism.
\item
A {\em morphism of left $G$-torsors\/} is just a morphism of the
underlying $G$-modules.  Likewise, we define right $G$-torsors,
$G$-bitorsors, and morphisms between them. We let:
$$
H^1(T,G)
$$
be the set of isomorphism classes of right $G$-torsors.
\item
A (left or right or bi-) $G$-torsor $(X,\mu_X)$ is said to be
{\em trivial}, if $\Gamma(T,X)\neq\emptyset$.
\end{enumerate}
\end{definition}

\begin{remark}\label{rem_standard-tors}
(i)\ \  In the situation of definition \ref{def_A-torsors},
notice that $H^1(T,G)$ always contains a distinguished element,
namely the class of the trivial $G$-torsor $(G,\mu_G)$.

(ii)\ \ Conversely, suppose that $(X,\mu_X)$ is a trivial left
$G$-torsor, and say that $\sigma\in\Gamma(T,X)$; then we have a
cartesian diagram :
$$\xymatrix{ G \ar[rr]^{\mu_\sigma} \ar[d] & &
X \ar[d]^{\one_X\times\sigma} \\
G\times X \ar[rr]^-{(\mu_X,p_X)} & & X\times X}
$$
which shows that $(X,\mu_X)$ is isomorphic to $(G,\mu_G)$ (and
likewise for right $G$-torsors).

(iii)\ \ Notice that every morphism $f:(X,\mu_X)\to(X',\mu_{X'})$ of
$G$-torsors is an isomorphism. Indeed, the assertion can be checked
locally on $T$ ({\em i.e.}, after pull-back by a covering morphism
$U\to 1_T$). Then we may assume that $X$ admits a global section
$\sigma\in\Gamma(T,X)$, in which case
$\sigma':=\sigma\circ f\in\Gamma(T,X')$. Then, arguing as in (ii),
we get a commutative diagram :
$$
\xymatrix{ & G \ar[ld]_{\mu_\sigma} \ar[rd]^{\mu_{\sigma'}} \\
X \ar[rr]^f & & X' }
$$
where both $\mu_\sigma$ and $\mu_{\sigma'}$ are isomorphisms,
and then the same holds for $f$.

(iv)\ \  The tensor product of a $G$-bitorsor and a left $G$-bitorsor
is a left $G$-torsor. Indeed the assertion can be checked locally on
$T$, so we are reduced to checking that the tensor product of a
trivial $G$-bitorsors and a trivial left $G$-torsor is the trivial
left $G$-torsor, which is obvious.

(v)\ \ Likewise, if $G_1\to G_2$ is any morphism of $T$-groups,
and $X$ is a left $G_1$-torsor, it is easily seen that the base
change $G_2\otimes_{G_1}X$ yields a left $G_2$-torsor (and the
same holds for right torsors and bitorsors). Hence the rule
$G\mapsto H^1(T,G)$ is a functor from the category of $T$-groups,
to the category of pointed sets. One can check that $H^1(T,G)$
is an essentially small set (see \cite[Chap.III, \S3.6.6.1]{Gi}).

(vi)\ \ Let $f:T\to S$ be a morphism of topoi; if $X$ is a left
$G$-torsor, the $f_*G$-module $f_*X$ is not necessarily a
$f_*G$-torsor, since we may not be able to find a covering morphism
$U\to 1_S$ such that $f_*X(U)\neq\emptyset$. On the other hand, if
$H$ a $S$-group and $Y$ a left $H$-torsor, then it is easily seen
that $f^*Y$ is a left $f^*H$-torsor.
\end{remark}

\sset\subsubsection{}\label{subsec_define-taus}
Let $f:T'\to T$ be a morphism of topoi, and $G$ a $T'$-monoid;
we define a $\sU$-presheaf $R^1f^\wedge_*G$ on $T$, by the rule :
$$
U\mapsto H^1(T'\!/\!f^*U,G_{|f^*U}).
$$
(More precisely, since this set is only essentially small,
we replace it by an isomorphic small set).
If $\phi:U\to V$ is any morphism in $T$, and $X$ is any
$G_{|f^*V}$-torsor, then $X\times_{f^*V}f^*U$ is a
$G_{|f^*U}$-torsor, whose isomorphism class depends only
on the isomorphism class of $X$; this defines the map
$R^1f^\wedge_*G(\phi)$, and it is clear that
$R^1f^\wedge_*G(\phi\circ\psi)=R^1f^\wedge_*G(\psi)\circ
R^1f^\wedge_*G(\phi)$, for any other morphism $\psi:W\to U$ in $T$.
Lastly, we denote by :
$$
R^1f_*G
$$
the sheaf on $(T,C_T)$ associated with the presheaf $R^1f^\wedge_*G$.
Notice that the object $R^1f_*G$ is {\em pointed}, {\em i.e.} it is
endowed with a natural global section :
$$
\tau_{f,G}:1_T\to R^1f_*G
$$
namely, the morphism associated with the morphism of presheaves
$1_T\to R^1f^\wedge_*G$ which, for every $U\in\Ob(T)$, singles out
the isomorphism class $\tau_{f,G}(U)\in R^1f^\wedge_*G(U)$ of the
trivial $G_{|f^*U}$-torsor.

\sset\subsubsection{}\label{subsec_Leray} Let $g:T''\to T'$ be
another morphism of topoi, and $G$ a $T''$-group. Notice that :
$$
f^{*\wedge}R^1g^\wedge_*G=R^1(f\circ g)^\wedge_*G
$$
hence the natural morphism (in $T^{\prime\wedge}$)
$R^1f^\wedge_*G\to R^1f_*G$ induces a morphism
$R^1(f\circ g)^\wedge_*G\to f^{*\wedge}R^1g_*G$ in $T^\wedge$,
which yields, after taking associated sheaves, a morphism in $T$ :
\set\begin{equation}\label{eq_Leray-I}
R^1(f\circ g)_*G\to f_*R^1g_*G.
\end{equation}
One sees easily that this is a {\em morphism of pointed objects\/}
of $T$, {\em i.e.} the image of the global section
$\tau_{f\circ g,G}$ under this map, is the global section
$f_*\tau_{g,G}$.

Next, suppose that $U\in\Ob(T)$ and $X$ is any $g_*G_{|f^*U}$-torsor
(on $T'/f^*U$); we may form the $g^*g_*G_{|g^*f^*U}$-torsor $g^*X$,
and then base change along the natural morphism $g^*g_*G\to G$, to
obtain the $G_{|g^*f^*U}$-torsor $G\otimes_{g^*g_*G}g^*X$.
This rule yields a map
$R^1f^\wedge_*(g_*G)\to R^1(f\circ g)^\wedge_*G$, and after taking
associated sheaves, a natural morphism of pointed objects :
\set\begin{equation}\label{eq_Leray-II}
R^1f_*(g_*G)\to R^1(f\circ g)_*G.
\end{equation}

\begin{remark}\label{rem_Leray-special}
As a special case, let $h:S'\to S$ be a morphism of topoi,
$H$ a $S'$-group. If we take $T'':=S'$, $T':=S$, $g:=h$ and
$f:S\to\Set$ the (essentially) unique morphism of topoi,
\eqref{eq_Leray-II} and \eqref{eq_Leray-I} boil down to maps
of pointed sets :
\set\begin{equation}\label{eq_special-Leray}
H^1(S,h_*H)\to H^1(S',H)\to\Gamma(S,R^1h_*H).
\end{equation}
\end{remark}

These considerations are summarized in the following :

\begin{theorem}\label{th_Leray-H^1}
In the situation of \eqref{subsec_Leray}, there exists a natural
{\em exact sequence of pointed objects} of\/ $T$ :
$$
1_T\to R^1f_*(g_*G)\to R^1(f\circ g)_*G\to f_*R^1g_*G.
$$
\end{theorem}
\begin{proof} The assertion means that \eqref{eq_Leray-II}
identifies $R^1f_*(g_*G)$ with the subobject :
$$
R^1(f\circ g)_*G\times_{f_*R^1g_*G}f_*\tau_{g,G}
$$
(briefly : the preimage of the trivial global section). 
We begin with the following :

\begin{claim}\label{cl_first-with-h}
In the situation of remark \ref{rem_Leray-special}, the
sequence of maps \eqref{eq_special-Leray} identifies the
pointed set $H^1(S,h_*H)$ with the preimage of the trivial
global section $\tau_{h,H}$ of $R^1h_*H$.
\end{claim}
\begin{pfclaim} Notice first that a global section $Y$ of
$R^1h_*^\wedge H$ maps to the trivial section $\tau_{h,H}$
of $R^1h_*H$ if and only if there exists a covering morphism
$U\to 1_S$ in $(S,C_S)$, such that $Y(h^*U)\neq\emptyset$.
Thus, let $X$ be a right $h_*H$-torsor; the image in $H^1(S',H)$
of its isomorphism class is the class of the $H$-torsor
$Y:=h^*X\otimes_{h^*h_*H}H$. The latter defines a global section
of $R^1h_*H$. However, by definition there exists a covering
morphism $U\to 1_S$ such that $X(U)\neq\emptyset$, hence also
$h^*X(h^*U)\neq\emptyset$, and therefore $Y(h^*U)\neq\emptyset$.
This shows that the image of $H^1(S,h_*H)$ lies in the preimage
of $\tau_{h,H}$.

Moreover, notice that $h_*Y(U)\neq\emptyset$, hence $h_*Y$ is a
$h_*H$-torsor. Now, let $\eps_H:h^*h_*H\to H$ (resp.
$\eta_{h_*H}:h_*H\to h_*h^*h_*H$) be the counit (resp. unit
of adjunction); we have a natural morphism
$\alpha:h^*X\to Y_{(\eps_H)}$ of $h^*h_*H$-modules, whence a
morphism :
$$
h_*\alpha:h_*h^*X\to h_*Y_{(h_*\eps_H)}
$$
of $h_*h^*h_*H$-modules. On the other hand, the unit of
adjunction $\eta_X:X\to h_*h^*X_{(\eta_{h_*H})}$ is a morphism
of $h_*H$-modules (remark \ref{rem_push-pull-mond}(v)). Since
$h_*\eps_H\circ\eta_{h_*H}=\one_{h_*H}$ (see \eqref{subsec_adj-pair}),
the composition $h_*\alpha\circ\eta_X$ is a morphism of
$h_*H$-modules, hence it is an isomorphism, by remark
\ref{rem_standard-tors}(iii). This implies that the first
map of \eqref{eq_special-Leray} is injective.

Conversely, suppose that the class of a $H$-torsor $X'$ gets
mapped to $\tau_{h,H}$; we need to show that the class of $X'$
lies in the image of $H^1(S,h_*H)$. However, the assumption
means that there exists a covering morphism $U\to 1_S$ such
that $X'(h^*U)\neq\emptyset$; by adjunction we deduce that
$h_*X'(U)\neq\emptyset$, hence $h_*X'$ is a $h_*H$-torsor.
In order to conclude, it suffices to show that the image in
$H^1(S',H)$ of the class of $h_*X'$ is the class of $X$.

Now, the counit of adjunction $h^*h_*X'\to X'$ is a morphism
of $h^*h_*H$-modules (remark \ref{rem_push-pull-mond}(v)); by
adjunction it induces a map $h^*h_*X'\otimes_{h^*h_*H}H\to H$
of $H$-torsors, which must be an isomorphism, according to
remark \ref{rem_standard-tors}(iii).
\end{pfclaim}

If we apply claim \ref{cl_first-with-h} with $S:=T'/f^*U$,
$S':=T''/(g^*f^*U)$ and $h:=g/(g^*f^*U)$, for $U$ ranging
over the objects of $T$, we deduce an exact sequence of
presheaves of pointed sets :
$$
1_T\to R^1f^\wedge_*(g_*H)\to R^1(f\circ g)^\wedge_*G\to
f^{*\wedge}R^1g_*G
$$
from which the theorem follows, after taking associated sheaves.
\end{proof}

\sset\subsubsection{}\label{subsec_coh_of_tors}
Let $f:T'\to T$ be a morphism of topoi, $U$ an object of $T$,
$G$ a $T'$-group, $p:X\to U$ a right $G_{|U}$-torsor.
Then, for every object $V$ of $T$ we have an induced
sequence of maps of sets
$$
X(f^*V)\xrightarrow{ p_* } U(f^*V)\xrightarrow{ \partial }
R^1f^\wedge_*G(V)
$$
where $p_*$ is deduced from $p$, and for every $\sigma\in U(f^*V)$,
we let $\partial(\sigma)$ be the isomorphism class of the right
$G_{|f^*V}$-torsor $(X\times_U f^*V\to f^*V)$. Clearly the image
of $p_*$ is precisely the preimage of (the isomorphism class of)
the trivial $G_{|\!f^*V}$-torsor. After taking associated sheaves,
we deduce a natural sequence of morphisms in $T$ :
\set\begin{equation}\label{eq_non-ab-coh}
f_*X \xrightarrow{ p } f_*U \xrightarrow{ \partial } R^1f_*G
\end{equation}
such that the preimage of the global section $\tau_{f,G}$ is
precisely the image in $\Gamma(T',U)$ of the set of global
sections of $X$.

\sset\subsubsection{}\label{subsec_topolo-on-sch}
A {\em ringed topos\/} is a pair $(T,\cO_T)$ consisting of
a topos $T$ and a (unital, associative) $T$-ring $\cO_T$,
called the {\em structure ring\/} of $T$. A morphism
$f:(T,\cO_T)\to(S,\cO_S)$ of ringed topoi is the datum of a
morphism of topoi $f:T\to S$ and a morphism of $T$-rings :
$$
f^\natural:f^*\cO_S\to\cO_T.
$$
We denote, as usual, by $\cO^\times_T\subset\cO_T$ the
subobject representing the invertible sections of $\cO_T$.
For every object $U$ of $T$, and every $s\in\cO_T(U)$,
let $D(s)\subset U$ be the subobject such that :
$$
\Hom_T(V,D(s)):=\{\phi\in U(V)~|~\phi^*s\in\cO^\times_T(V)\}.
$$
We say that $(T,\cO_T)$ is {\em locally ringed}, if 
$D(0)=\emptyset_T$ (the initial object of $T$), and moreover
$$
D(s)\cup D(1-s)=U
\qquad
\text{for every $U\in\Ob(T)$, and every $s\in\cO_T(f)$.}
$$
A morphism of locally ringed topoi $f:(T,\cO_T)\to(S,\cO_S)$
is a morphism of ringed topoi such that
$$
f^*D(s)=D(f^\natural(U)(f^*s))
\qquad
\text{for every $U\in\Ob(S)$ and every $s\in\cO_S(U)$.}
$$
If $T$ has enough points, then $(T,\cO_T)$ is locally ringed
if and only if the stalks $\cO_{T,\xi}$ of the structure ring
at all the points $\xi$ of $T$ are local rings. Likewise, a
morphism $f:(T,\cO_{\!T})\to(S,\cO_{\!S})$ of ringed topoi is
locally ringed if and only if, for every $T$-point $\xi$,
the induced map $\cO_{\!S,f(\xi)}\to\cO_{\!T,\xi}$ is a
local ring homomorphism.

\sset\subsubsection{}\label{subsec_etale-and-Zariski}
In the rest of this section we present a few results concerning the
special case of topologies on a scheme. Hence, for any scheme $X$,
we shall denote by $X_\et$ (resp. by $X_\Zar$) the small {\'e}tale
(resp. the small Zariski) site on $X$. It is clear that $X_\Zar$
is a small site, and it is not hard to show that $X_\et$ is a
$\sU$-site (\cite[Exp.VII, \S1.7]{SGA4-2}).
The inclusion of underlying categories :
$$
u_X:X_\Zar\to X_\et
$$
is a continuous functor (see definition \ref{def_dir-img-site}(i))
commuting with finite limits, whence a morphism of topoi :
$$
\tilde u_X:=(\tilde u{}^*_{\!X},\tilde u_{X*}):X^\sim_\et\to
X^\sim_\Zar.
$$
such that the diagram of functors :
$$
\xymatrix{ X_\Zar \ar[r]^-{u_X} \ar[d] & X_\et \ar[d] \\
           X^\sim_\Zar \ar[r]^-{\tilde u{}^*_{\!X}} & X^\sim_\et
}$$ commutes, where the vertical arrows are the Yoneda embeddings
(lemma \ref{lem_cont-funct-site}).

The topoi $X^\sim_\Zar$ and $X^\sim_\et$ are locally ringed in a
natural way, and by faithfully flat descent, we see easily that
$\tilde u_{X*}\cO_{\!X_\et}=\cO_{\!X_\Zar}$. By inspection,
$\tilde u_X$ is a morphism of locally ringed topoi.

\sset\subsubsection{}\label{subsec_big-site}
For any ring $R$ (of our fixed universe $\sU$), denote by
$\Sch/R$ the category of $R$-schemes, and by $\Sch/R_\Zar$ (resp.
$\Sch/R_\et$) the big Zariski (resp. \'etale) site on $\Sch/R$.
For $R=\Z$, we shall usually just write $\Sch_\Zar$ and $\Sch_\et$
for these sites.
The morphisms $u_X$ of \eqref{subsec_topolo-on-sch} are actually
restrictions of a single morphism of sites :
$$
u:\Sch_\Zar\to\Sch_\et
$$
which, for every universe $\sV$ such that $\sU\in\sV$,
induces a morphism of $\sV$-topoi :
$$
\tilde u_\sV:(\Sch_\et)_\sV^\sim\to(\Sch_\Zar)_\sV^\sim.
$$

\sset\subsubsection{}\label{subsec_geom-pts}
Let $X$ be scheme; a {\em geometric point\/} of $X$ is a
morphism of schemes $\xi:\Spec\,\kappa\to X$, where $\kappa$
is an arbitrary separably closed field. Notice that both the
Zariski and \'etale topoi of $\Spec\,\kappa$ are equivalent
to the category $\Set$, so $\xi$ induces a  topos-theoretic
point $\xi_\et^\sim:(\Spec\,\kappa)^\sim_\et\to X^\sim_\et$
of $X^\sim_\et$ (and likewise for $X^\sim_\Zar$). A basic
feature of both the Zariski and \'etale topologies, is that
every point of $X^\sim_\Zar$ and $X^\sim_\et$ arise in this
way.

More precisely, we say that two geometric points $\xi$ and
$\xi'$ of $X$ are {\em equivalent}, if there exists a third
such point $\xi''$ which factors through both $\xi$ and $\xi'$.
It is easily seen that this is an equivalence relation on the
set of geometric points of $X$, and two topos-theoretic points
$\xi^\sim_\et$ and $\xi^{\prime\sim}_\et$ are isomorphic if and
only if the same holds for the points $\xi^\sim_\Zar$ and
$\xi^{\prime\sim}_\Zar$, if and only if $\xi$ is equivalent
to $\xi'$.

\begin{definition}\label{def_strict-loc}
Let $X$ be a scheme, $x$ a point of $X$, and
$\bar x:\Spec\,\kappa\to X$ a geometric point.

(i)\ \
We let $\kappa(x)$ be the residue field of the local ring
$\cO_{\!X,x}$, and set
$$
|x|:=\Spec\,\kappa(x)
\qquad
\kappa(\bar x):=\kappa
\qquad
|\bar x|:=\Spec\,\kappa(\bar x)
$$
If $\{x\}\subset X$ is the image of $\bar x$, we say that $\bar x$
is {\em localized at $x$}, and that $x$ is the {\em support\/} of
$\bar x$. We say that $\bar x$ is {\em strict\/} if $\kappa(\bar x)$
is a separable closure of $\kappa(x)$.

(ii)\ \
We associate with $\bar x$ a strict geometric point $\bar x{}^\rmst$,
as follows. Let $\kappa(\bar x{}^\rmst)$ be the separable closure of
$\kappa(x)$ inside $\kappa(\bar x)$; the inclusion
$\kappa(\bar x{}^\rmst)\subset\kappa(\bar x)$ defines a morphism
of schemes :
\set\begin{equation}\label{eq_strictify}
|\bar x|\to|\bar x{}^\rmst|:=\Spec\,\kappa(\bar x{}^\rmst)
\end{equation}
and $\bar x$ is the composition of \eqref{eq_strictify} and
a unique strict geometric point localized at $x$
$$
\bar x{}^\rmst:|\bar x{}^\rmst|\to X.
$$

(iii)\ \
The {\em localization of $X$ at $x$} is the local scheme
$$
X(x):=\Spec\,\cO_{\!X,x}.
$$
The {\em strict henselization of $X$ at $\bar x$} is the
strictly local scheme
$$
X(\bar x):=\Spec\,\cO_{\!X,\bar x}
$$
where $\cO_{\!X,\bar x}$ denotes the strict henselization of
$\cO_{\!X,x}$ relative to the geometric point $\bar x$
(\cite[Ch.IV, D{\'e}f.18.8.7]{EGA4}) (recall that a local ring
is called {\em strictly local}, if it is henselian with separably
closed residue field; a scheme is called {\em strictly local}, if
it is the spectrum of a strictly local ring: see
\cite[Ch.IV, D{\'e}f.18.8.2]{EGA4}). By definition, the
geometric point $\bar x$ lifts to a unique geometric point
of $X(\bar x)$, which shall be denoted again by $\bar x$.

(iv)\ \
Moreover, we shall denote by
$$
i_x:X(x)\to X \qquad i_{\bar x}:X(\bar x)\to X(x)
$$
the natural morphisms of schemes, and if $\cF$ is any
sheaf on $X_\Zar$ (resp. $X_\et$), we let
$$
\cF(x):=i_x^*\cF
\qquad
\cF(\bar x):=i_{\bar x}^*\cF(x)
$$
so $\cF(x)$ is a sheaf on $X(x)_\Zar$ (resp. $X(x)_\et$)
and $\cF(\bar x)$ is a sheaf on $X(\bar x)_\Zar$ (resp.
$X(\bar x)_\et$).

(v)\ \
If $f:Y\to X$ is any morphism of schemes, we let
$$
f^{-1}(x):=Y\times_X|x|
\quad
f^{-1}(\bar x):=Y\times_X|\bar x|
\qquad
Y(x):=Y\times_XX(x)
\quad
Y(\bar x):=Y\times_XX(\bar x).
$$
Also, if $\bar y$ is any geometric point of $Y$, we define
$f(\bar y)$ as the geometric point $f\circ\bar y$ of $X$, and
we call $f(\bar y)^\rmst$ the {\em strict image\/} of $\bar y$
in $X$. Notice that the natural identification
$|\bar y|=|f(\bar y)|$ induces a morphism of schemes :
$$
|\bar y{}^\rmst|\to|f(\bar y)^\rmst|.
$$
\end{definition}

\sset\subsubsection{}\label{subsec_tau-notate}
Many discussions concerning the Zariski or \'etale site of a
scheme, only make appeal to general properties of these two
topologies, and therefore apply indifferently to either of
them, with only minor verbal changes. For this reason, to
avoid tiresome repetitions, the following notational device
is often useful. Namely, instead of referring each time to
$X_\Zar$ and $X_\et$ in the course of an argument, we shall
write just $X_\tau$, with the convention that $\tau\in\{\Zar,\et\}$
has been chosen arbitrarily at the beginning of the discussion.
In the same manner, a $\tau$-point of $X$ will mean a point of the
topos $X^\sim_\tau$, and a $\tau$-open subset of $X$ will be any
object of the site $X_\tau$. With this convention, a $\Zar$-point
is a usual point of $X$, whereas an $\et$-point shall be a geometric
point. Likewise, if $\xi$ is a given $\tau$-point of $X$, the
localization $X(\xi)$ makes sense for both topologies : if
$\tau=\et$, then $X(\xi)$ is the strict henselization as in
definition \ref{def_strict-loc}(ii); if $\tau=\Zar$, then $X(\xi)$
is the usual localization of $X$ at the (Zariski) point $\xi$.
If $\tau=\et$, the support of $\xi$ is given by definition
\ref{def_strict-loc}(i); if $\tau=\Zar$, then the support of $\xi$
is just $\xi$ itself (and correspondingly, in this case $\xi$
is localized at $\xi$). Furthermore, $\cO_{\!X,\xi}$ is a local
ring if $\tau=\Zar$, and it is a strictly local ring, in case
$\tau=\et$.

\sset\subsubsection{}\label{subsec_strict-loc-of-schs}
Let $f:X\to Y$ be a morphism of schemes, $\bar x$ a geometric
point of $X$, and set $\bar y:=f(\bar x)$. The natural morphism
$f_x:X(x)\to Y(y)$ induces a unique local morphism of strictly
local schemes
$$
f_{\bar x}:X(\bar x)\to Y(\bar y)
$$
(\cite[Ch.IV, Prop.18.8.8(ii)]{EGA4}) that fits in a commutative diagram :
$$
\xymatrix{
|\bar x| \ar[r]^-{\bar x} \ddouble &
X(\bar x) \ar[d]_{f_{\bar x}} \ar[r]^{i_{\bar x}} &
X(x) \ar[d]^{f_x} \ar[r]^-{i_x} & X \ar[d]^f \\
|\bar y| \ar[r]^-{\bar y} & Y(\bar y) \ar[r]^{i_{\bar y}} & Y(y)
\ar[r]^-{i_y} & Y. }
$$
Let now $\cF$ be any sheaf on $Y_\et$; there follows a natural
isomorphism :
$$
f^*_{\bar x}\cF(\bar y)\isom(f^*\cF)(\bar x).
$$
Notice also the natural bijections :
\set\begin{equation}\label{eq_maps}
\begin{aligned}
\cF_{\bar y} & \isom\Gamma(Y(\bar y),\cF(\bar y))\isom
\Gamma(|\bar y|,\bar y{}^*\cF(\bar y)) \\
f^*\cF_{\bar x} & \isom\Gamma(X(\bar x),f^*\cF(\bar x))\isom
\Gamma(|\bar x|,\bar x{}^*f^*\cF(\bar x))
\end{aligned}
\end{equation}
which induce a natural identification :
\set\begin{equation}\label{eq_madness}
\cF_{\bar y}\isom f^*\cF_{\bar x}
\quad :\quad
\sigma\mapsto f^*_{\bar x}(\sigma).
\end{equation}

\sset\subsubsection{}\label{subsec_strict-special}
Let $X$ be a scheme, $x,x'\in X$ any two points, such that
$x$ is a specialization of $x'$. Choose a geometric point
$\bar x$ localized at $x$.
The localization map $\cO_{\!X,x}\to\cO_{\!X,x'}$ induces
a natural {\em specialization morphism\/} of $X$-schemes :
$$
X(x')\to X(x).
$$
Set $W:=X(\bar x)\times_{X(x)}X(x')$. The natural map
$g:W\to X(x')$ is faithfully flat, and is the limit of a
cofiltered system of \'etale morphisms; hence we may
find $w\in W$ lying over the closed point of $X(x')$,
and the induced map $\kappa(x')\to\kappa(w)$ is algebraic
and separable. Choose also a geometric point $\bar w$
of $W$ localized at $w$, and set $\bar x{}':=g(\bar w)$.
Then $g$ induces an isomorphism
$g_{\bar w}:W(\bar w)\isom X(\bar x{}')$, whence a unique
morphism
\set\begin{equation}\label{eq_speci-map}
X(\bar x{}')\to X(\bar x)
\end{equation}
which makes commute the diagram :
$$
\xymatrix{
W(\bar w) \ar[r]^-{g_{\bar w}} \ar[d]_{i_{\bar w}} &
X(\bar x{}') \ar[d] \ar[r]^-{i_{\bar x{}'}} & X(x') \ar[d] \\
W(w) \ar[r] & X(\bar x) \ar[r]^{i_{\bar x}} & X(x)
}$$
where the left bottom arrow is the natural projection, and
the right-most vertical arrow is the specialization map. 
In this situation, we say that $\bar x$ is a
{\em specialization\/} of $\bar x{}'$ (and that $\bar x{}'$
is a {\em generization\/} of $\bar x$), and we call
\eqref{eq_speci-map} a {\em strict specialization morphism}.
Combining with \eqref{eq_maps}, we obtain the
{\em strict specialization map induced by} \eqref{eq_speci-map}
\set\begin{equation}\label{eq_spec-map-induc}
\cG_{\bar x}\to\cG_{\bar x{}'}
\end{equation}
for every sheaf $\cG$ on $X_\et$. 

\begin{remark}\label{rem_down-stri-spec}
(i)\ \
In the situation of \eqref{subsec_strict-loc-of-schs},
suppose that $\cG=f^*\cF$ for a sheaf $\cF$ on $Y_\et$.
Then \eqref{eq_spec-map-induc} is a map
$\cF_{f(\bar x)}\to\cF_{f(\bar x{}')}$. By inspecting the
definition, it is easily seen that the latter agrees with
the strict specialization map for $\cF$ induced by a
unique strict specialization morphism
$Y(f(\bar x{}'))\to Y(f(\bar x))$.

(ii)\ \
Notice that \eqref{eq_speci-map} and \eqref{eq_spec-map-induc}
depend not only on the choice of $w$ (which may not be unique,
when $X(x)$ is not unibranch) but also on the geometric point
$\bar w$. Indeed, the group of automorphisms of the $X(x')$-scheme
$X(\bar x{}')$ is naturally isomorphic to the Galois group
$\Gal(\kappa(x')^\mathrm{s}/\kappa(x'))$
(\cite[Ch.IV, (18.8.8.1)]{EGA4}).
\end{remark}

\begin{lemma}\label{lem_Hilbert90}
Let $X$ be a scheme, $\cF$ a sheaf on $X_\et$. We have :
\begin{enumerate}
\item
The counit of the adjunction $\eps_\cF:\tilde
u{}^*_{\!X}\circ\tilde u_{X*}\cF\to\cF$ is a monomorphism.
\item
Suppose there exist a sheaf\/ $\cG$ on $X_\Zar$, and an epimorphism
$f:\tilde u{}^*\cG\to\cF$ (resp. a monomorphism $f:\cF\to\tilde u{}^*\cG$).
Then $\eps_\cF$ is an isomorphism.
\item
The functor $\tilde u{}^*_{\!X}$ is fully faithful.
\item
{\em (Hilbert 90)}\ \ $R^1\tilde u_{X*}\cO^\times_{\!X_\et}=1_{X_\Zar}$.
\end{enumerate}
\end{lemma}
\begin{proof}(i): The assertion can be checked on the stalks. Hence,
let $\xi$ be any geometric point of $X$; we have to show that the
natural map $(\tilde u_{X*}\cF)_\xi\to\cF_\xi$ is injective. To this
aim, say that $s,s'\in(\tilde u_{X*}\cF)_\xi$, and suppose that the
image of $s$ in $\cF_\xi$ agrees with the image of $s'$; we may find
an open neighborhood $U$ of $\xi$ in $X_\Zar$, such that $s$ and
$s'$ lie in the image of $\cF(U)$, and by assumption, there exists
an \'etale morphism $f:V\to U$ such that the images of $s$ and $s'$
coincide in $\cF(V)$. However, $f(V)\subset U$ is an open subset
(\cite[Ch.IV, Th.2.4.6]{EGAIV-2}), and the induced map $V\to f(V)$
is a covering morphism in $X_\et$; it follows that the images of $s$
and $s'$ agree already in $\cF(f(V))$, therefore also in $(\tilde
u_{X*}\cF)_\xi$.

(iii): According to proposition \ref{prop_fullfaith-adjts}(iii),
it suffices to show that the unit of the adjunction
$\eta_\cG:\cG\to\tilde u_{X*}\circ\tilde u{}^*_{\!X}\cG$
is an isomorphism, for every $\cG\in\Ob(X^\sim_\Zar)$. However,
we have morphisms :
$$
\tilde u{}^*_X\cG\xrightarrow{\ \tilde u{}^*_X(\eta_\cG) }
\tilde u{}^*_X\circ\tilde u_{X*}\circ\tilde u{}^*_{\!X}\cG 
\xrightarrow{\ \eps_{\tilde u{}^*_X\cG} }\tilde u{}^*_X\cG.
$$
whose composition is the identity of $\tilde u{}^*_X\cG$
(see \eqref{subsec_adj-pair}); also, (i) says that
$\eps_{\tilde u{}^*_X\cG}$ is a monomorphism, and then
it follows formally that it is actually an isomorphism
({\em e.g.} from the dual of \cite[Prop.1.9.3]{Bor}).
Hence the same holds for $\tilde u{}^*_X(\eta_\cG)$, and by
considering the stalks of the latter, we conclude that also
$\eta_\cG$ is an isomorphism, as required.

(ii): Suppose first that $f:\tilde u{}^*\cG\to\cF$ is an
epimorphism. We have just seen that $\eta_\cG$ is an
isomorphism, therefore we have a morphism
$\tilde u{}^*\circ\tilde u_*f:
\tilde u{}^*\cG\to\tilde u{}^*\circ\tilde u_*\cF$ whose
composition with $\eps_\cF$ is $f$; especially, $\eps_\cF$
is an epimorphism, so the assertion follows from (i).

In the case of a monomorphism $f:\cF\to\tilde u{}^*\cG$,
set $\cH:=\tilde u^*\cG\amalg_\cF\tilde u^*\cG$;
we may represent $f$ as the equalizer of the two natural
maps $j_1,j_2:\tilde u{}^*\cG\to\cH$. However, the natural
morphism $\tilde u{}^*(\cG\amalg\cG)\to\cH$ is an epimorphism,
hence the counit $\eps_\cH$ is an isomorphism, by the previous
case. Then (iii) implies that $j_i=\tilde u{}^*j'_i$ for
morphisms $j'_i:\cG\to\tilde u_*\cH$ ($i=1,2$). Let $\cF'$
be the equalizer of $j'_1$ and $j'_2$; then
$\tilde u{}^*\cF'\simeq\cF$, and since we have already seen
that the unit of adjunction is an isomorphism, the assertion
follows from the triangular identitities of \eqref{subsec_adj-pair}.

(iv): The assertion can be checked on the stalks. To ease notation,
set $\cF:=R^1\tilde u_{X*}\cO^\times_{\!X_\et}$. Let $\xi$ be any
geometric point of $X$, and say that $s\in\cF_\xi$; pick a (Zariski)
open neighborhood $U\subset X$ of $\xi$ such that $s$ lies in the image
of $\cF(U)$. We may then find a Zariski open covering $(U_\lambda\to
U~|~\lambda\in\Lambda)$ of $U$, such that the image of $s$ in
$\cF(U_i)$ is represented by a $\cO^\times_{U_\lambda,\et}$-torsor
on $U_{\lambda,\et}$, for every $\lambda\in\Lambda$. After replacing
$U$ by any $U_\lambda$ containing the support of $\xi$, we may
assume that $s$ is the image of the isomorphism class of some
$\cO_{\!U_\et}^\times$-torsor $X_\et$ on $U_\et$. By faithfully flat
descent, there exist a $\cO_{\!U_\Zar}^\times$-torsor $X$ on
$U_\Zar$, and an isomorphism of $\cO_{\!U_\et}^\times$-torsors :
$$
X_\et\isom\cO_{\!U_\et}^\times\otimes_{\tilde
u{}^*_U\cO_{\!U_\Zar}^\times}\tilde u{}^*_UX.
$$
However, after replacing $U$ by a smaller open neighborhood of
$\xi$, we may suppose that $X(U)\neq\emptyset$, therefore
$X_\et(U)\neq\emptyset$ as well, {\em i.e.} $s$ is the image of the
trivial section of $\cF(U)$.
\end{proof}

\section{Stacks}
Let $(\cC,J)$ be any site, and consider the rule that assigns
to every $U\in\Ob(\cC)$ the category $(U,J_U)^\sim$ of sheaves
on $\cC/U$, for the topology $J_U$ induced by $J$ via the source
functor $\cC/U\to\cC$. Every morphism $g:U'\to U$ of $\cC$ induces
a pull-back functor $\tilde g{}^*:(U,J_U)^\sim\to(U',J_{U'})^\sim$,
but if $g':U''\to U'$ is another such morphism, the composition
$\tilde g'{}^*\circ\tilde g{}^*$ is not usually equal, but only
{\em isomorphic} to the pull-back functor induced by $g\circ g'$.
Moreover, given a covering family $(U_i\to U~|~i\in I)$, we can
describe the category $(U,J_U)^\sim$ only {\em up to equivalence}
in terms of the corresponding categories on the $U_i$ and their
intersections; this ``gluing up to equivalence'' is a
$2$-categorical manipulation of descent data relative to the
given covering, so it involves not only double intersections
$U_i\times U_j$, but also triple intersections. Summing up, the
rule $U\mapsto(U,J_U)^\sim$ does not quite define a sheaf of
categories on $C$, but rather a kind of $2$-categorical analogue
of the latter : it is the first example of a structure that is
often encountered in algebraic geometry, and whose systematic study
is the subject of the {\em theory of stacks} developed in this chapter.

\subsection{Prestacks and stacks on a site}
As usual, $\sU$ will denote a chosen universe, and small
will be synonymous with $\sU$-small, throughout this section.

\begin{definition} Let $C:=(\cC,J)$ be a site, $\phi:\cA\to\cC$
a fibration, and $i\leq 2$ an integer. We say that $\phi$ is
an {\em $i$-separated prestack} if for every $X\in\Ob(\cC)$,
every $\cS\in J(X)$ is a sieve of $\phi$-$i$-descent
(see definition \ref{def_descent-fibred}). A $2$-separated
prestack is called a {\em stack}. For every universe $\sV$,
we denote by
$$
\sV\tdu i\tdu\PreStack(C)
\qquad
\text{(resp.\ \ $\sV\tdu\Stack(C)$\ )}
$$
the full $2$-subcategory of $\sV\tdu\Fib(\cC)$ whose objects
are the $i$-separated prestacks (resp. the stacks) on $C$.
Hence $\sV\tdu(-1)\tdu\PreStack(C)=\sV\tdu\Fib(\cC)$ and
$\sV\tdu 2\tdu\PreStack(C)=\sV\tdu\Stack(C)$. As usual,
if $\sV=\sU$ we often drop the mention of the universe
from this notation.
\end{definition}

\begin{example}\label{ex_sheaves-as-stacks}
Let $(\cC,J)$ be any site, $F$ any presheaf on $\cC$, and
$\ss_F:\cFib(F)\to\cC$ the fibration attached to $F$, as
in \eqref{subsec_fibred-cats-II}.
Let also $X\in\Ob(\cC)$ be any object, and $\cS$ any sieve
of $\cC/X$; by \eqref{subsec_fibred-cats-II} and remark
\ref{rem_sieves-and-sub}(ii) we have a commutative
diagram
$$
\xymatrix{ FX=\Hom_{\cC^\wedge}(h_X,F) \ar[r] \ar[d] &
\Hom_{\cC^\wedge}(h_\cS,F) \ar[d] \\
\cFib(F)(X) \ar[r] & \sCart_\cC(\cS,\cFib(F))
}$$
whose vertical arrows are bijections, and whose top (resp.
bottom) horizontal arrow is induced by the inclusion
$h_\cS\to h_X$ (resp. $\cS\subset\cC/X$). We conclude that
$\cFib(F)$ is a $0$-separated prestack; by the same token,
$\cFib(F)$ is a $1$-separated prestack (resp. a stack) if
and only if $F$ is a separated presheaf (resp. a sheaf).
Thus, we get a commutative diagram of functors
\set\begin{equation}\label{eq_sheaves-as-stacks}
{\diagram
(\cC,J)^\sim \ar[r] \ar[d] & \Stack(\cC,J) \ar[d] \\
\cC^\wedge \ar[r]^-{\cFib} & 0\tdu\PreStack(\cC,J)
\enddiagram}
\end{equation}
(whose vertical arrows are the inclusion functors) which
allows us to regard $(\cC,J)^\sim$ as a full subcategory
of $\Stack(\cC,J)$.
\end{example}

\sset\subsubsection{}\label{subsec_presheaf-of-hom-of-carts}
We can characterize $0$-separated and $1$-separated prestacks on
a small site $C:=(\cC,J)$ by means of the following construction.
Let $\cA\to\cC$ be any fibration with small fibres. To every
$X\in\Ob(\cC)$ and every pair of cartesian sections
$\sigma,\sigma'\in\cA(X)$ we attach the presheaf
$$
\cCart(\sigma,\sigma'):(\cC/X)^o\to\Set
\qquad
(Y\xrightarrow{f}X)\mapsto\sCart_\cC(\sigma\circ f_*,\sigma'\circ f_*)
$$
which assigns to every morphism
$(Z\xrightarrow{g}X)\xrightarrow{h/X}(Y\xrightarrow{f}X)$ of $\cC/X$
the map
$$
\cCart(\sigma,\sigma')(f)\to\cCart(\sigma,\sigma')(g)
\qquad
\beta\mapsto\beta*h_*
$$
(with $f_*:\cC/Y\to\cC/X$ and $h_*:\cC/Z\to\cC/Y$ as
in \eqref{eq_push-for}). Notice that :
\set\begin{equation}\label{eq_triviality}
(f_*)^\wedge\cCart(\sigma,\sigma')=
\cCart(\sigma\circ f_*,\sigma'\circ f_*)
\qquad
\text{for every $(Y\xrightarrow{f}X)\in\Ob(\cC/X)$}.
\end{equation}
Moreover, for every $\sigma,\sigma',\sigma''\in\cA(X)$,
the system of composition maps
$$
\sCart_\cC(\sigma\circ f_*,\sigma'\circ f_*)\times
\sCart_\cC(\sigma'\circ f_*,\sigma''\circ f_*)\to
\sCart_\cC(\sigma\circ f_*,\sigma''\circ f_*)
\qquad
(\beta,\beta')\mapsto\beta'\odot\beta
$$
clearly defines a morphism of presheaves on $\cC/X$ :
\set\begin{equation}\label{eq_surfing-Doc}
\cCart(\sigma,\sigma')\times\cCart(\sigma',\sigma'')
\to\cCart(\sigma,\sigma'').
\end{equation}
Furthermore, every pair of isomorphisms of cartesian
sections : $\mu:\sigma\isom\tau$, $\mu':\sigma'\isom\tau'$
induces an isomorphism of presheaves
\set\begin{equation}\label{eq_pariah}
\cCart(\sigma,\sigma')\isom\cCart(\tau,\tau')
\qquad
(\beta:\sigma\circ f_*\Rightarrow\sigma'\circ f_*)\mapsto
(\mu'*f_*)\odot\beta\odot(\mu^{-1}*f_*).
\end{equation}
Lastly, every cartesian functor $F:\cA\to\cA'$ of
$\cC$-fibrations yields a morphism of presheaves :
\set\begin{equation}\label{eq_F-pariah}
\cCart(\sigma,\sigma')\to\cCart(F\circ\sigma,F\circ\sigma')
\qquad
\beta\mapsto F*\beta.
\end{equation}

\begin{lemma}\label{lem_crit-0-1_separation}
With the notation of \eqref{subsec_presheaf-of-hom-of-carts},
the following conditions are equivalent :
\begin{enumerate}
\alphaenu
\item
The fibration $\cA\to\cC$ is $0$-separated (resp. $1$-separated).
\item
For every $X\in\Ob(\cC)$ and every pair of cartesian sections
$\sigma,\sigma'\in\cA(X)$, the presheaf\/ $\cCart(\sigma,\sigma')$
is separated (resp. is a sheaf) on the site $C/X$ (notation of
\eqref{sec_Localization-topoi}).
  \end{enumerate}
\end{lemma}
\begin{proof} Let $X\in\Ob(\cC)$ be any object, and $\cS\subset\cC/X$
any covering sieve; choose a generating family $(f_i:Y_i\to X~|~i\in I)$
for $\cS$, and for every $i,j\in I$ set
$\cC/Y_{ij}:=\cC/Y_i\times_{\cC/X}\cC/Y_j$. There follows a
commutative diagram :
$$
\xymatrix{ \cA(X) \ar[rr]^-{\sCart_\cC(j,\cA)}
\ar[rrd]_\rho & & \sCart_\cC(\cS,\cA) \ar[d]^{\eps^*} \\
& & \prod_{i\in I}\cA(Y_i)
\ar@<.5ex>[r]^-{\partial^*_0} \ar@<-.5ex>[r]_-{\partial^*_1} &
\prod_{i,j\in I}\sCart_\cC(\cC/Y_{ij},\cA)
}$$
where $j:\cS\to\cC/X$ is the inclusion functor, $\rho$ is the
unique functor whose composition with the projection onto the
factor $\cA(Y_i)$ agrees with $\sCart_\cC(f_{i*},\cA)$,
for every $i\in I$, and the functors $\eps^*$, $\partial^*_0$ and
$\partial^*_1$ are as in remark \ref{rem_fin-prod-reps}. Clearly
$\eps^*$ is faithful; it follows that $\rho$ is faithful if and
only if the same holds for $\sCart_\cC(j,\cA)$. On the other hand,
the family $((f_i/X):f_i\to\one_X~|~i\in I)$ of morphisms of
$\cC/X$ covers $\one_X$, according to \eqref{sec_Localization-topoi}.
Taking into account \eqref{eq_triviality}, we conclude easily
that $\cA$ is $0$-separated if and only if all the presheaves
$\cCart(\sigma,\sigma')$ are separated. Next, by remark
\ref{rem_fin-prod-reps} the functor $\eps^*$ induces an
isomorphism of $\sCart_\cC(\cS,\cA)$ onto the equalizer $E$
of $\partial^*_0$ and $\partial^*_1$, and clearly $\rho$
factors through a unique functor $\rho':\cA(X)\to E$. It follows
that $\rho'$ is fully faithful if and only if the same holds for
$\sCart_\cC(j,\cA)$. Again, this means that $\cA$ is $1$-separated
if and only if all the presheaves $\cCart(\sigma,\sigma')$
are sheaves.
\end{proof}

\begin{theorem}\label{th_be-positive}
Let $C:=(\cC,J)$ be any small site. The inclusion strict
pseudo-functor $1\tdu\PreStack(C)\to\Fib(\cC)$ admits a
strict and strong left $2$-adjoint pseudo-functor :
$$
\Fib(\cC)\to1\tdu\PreStack(C)
\qquad
\cA\mapsto\cA^\sep.
$$
\end{theorem}
\begin{proof} Let $\cA\to\cC$ be any fibration with small
fibre categories; pick a unital cleavage $\blambda$ for
$\cA$, and let $\sc$ be its associated unital pseudo-functor.
For every $X\in\Ob(\cC)$ consider the functor
$$
\beta^\blambda_X:\cA_X\to\cA(X)
\qquad
A\mapsto\beta^\blambda_{X,A}
$$
defined as in the proof of theorem \ref{th_split-fibration}.
For every $A,A'\in\Ob(\cA_X)$, consider as well the map
$$
\omega_{A,A'}:\Hom_{\cA_X}(A,A')\to
\cCart(\beta^\blambda_{X,A},\beta^\blambda_{X,A'})(\one_X)
\qquad
(f:A\to A')\mapsto\beta^\blambda_{X,f}.
$$
Let $\cH_{A,A'}$ be the sheaf on $C/X$ associated with the
presheaf $\cCart(\beta^\blambda_{X,A},\beta^\blambda_{X,A'})$,
and denote
$$
\tilde\omega_{A,A'}:\Hom_{\cA_X}(A,A')\to H_{A,A'}:=\cH_{A,A'}(\one_X)
$$
the composition of $\omega_{A,A'}$ with the natural map
$\cCart(\beta^\blambda_{X,A},\beta^\blambda_{X,A'})(\one_X)\to
H_{A,A'}$. For every $A,A',A''\in\Ob(\cA_X)$, we have morphisms
of presheaves as in \eqref{eq_surfing-Doc}
$$
\cCart(\beta^\blambda_{X,A},\beta^\blambda_{X,A'})\times
\cCart(\beta^\blambda_{X,A'},\beta^\blambda_{X,A''})\to
\cCart(\beta^\blambda_{X,A},\beta^\blambda_{X,A''})
$$
whence a morphism of sheaves
$\cH_{A,A'}\times\cH_{A',A''}\to\cH_{A,A''}$, which induces a
composition map :
$$
H_{A,A'}\times H_{A',A''}\to H_{A,A''}.
$$
It is easily seen that this composition law is
associative, for every $A,A',A'',A'''\in\Ob(\cA_X)$,
and $\tilde\omega_{A,A}(\one_A)$ is neutral for left
and right composition. Hence, we have a category $\cA_X^\sep$
whose set of objects is $\Ob(\cA_X)$, and such that
$\Hom_{\cA_X^\sep}(A,A'):=H_{A,A'}$ for every $A,A'\in\Ob(\cA_X)$,
with the foregoing composition law. Furthermore, the
system of maps $\tilde\omega_{\bullet\bullet}$ yields a functor
$$
\tilde\omega^X:\cA_X\to\cA^\sep_X
$$
which is the identity map on objects. Next, notice that
for every morphism $g:Y\to X$ in $\cC$ and every
$A\in\Ob(\cA_X)$ we have
$\beta^\blambda_{X,A}(g)=\beta^\blambda_{Y,\sc_gA}(\one_Y)=\sc_gA$;
it follows that there exists a unique isomorphism of
cartesian sections
\set\begin{equation}\label{eq_fart}
\beta^\blambda_{X,A}\circ g_*\isom\beta^\blambda_{Y,\sc_gA}
\qquad\text{such that}\qquad
\one_Y\mapsto\one_{\sc_gA}.
\end{equation}
By virtue of \eqref{eq_triviality} and \eqref{eq_pariah}, the
isomorphisms \eqref{eq_fart} induce an isomorphism of presheaves
$$
\mu^g_{(A,A')}:g^\wedge_*\cCart(\beta^\blambda_{X,A},\beta^\blambda_{X,A'})
\isom\cCart(\beta^\blambda_{Y,\sc_gA},\beta^\blambda_{Y,\sc_gA'})
\qquad
\text{for every $A,A'\in\Ob(\cA_X)$}.
$$
Since $g_*$ is continuous and cocontinuous for the topologies
of $C/X$ and $C/Y$ (remark \ref{rem_continue-local}(i)), combining
with lemma \ref{lem_improve}(ii) we deduce a natural isomorphism
of sheaves on $C/Y$ :
$$
\tilde\mu^g_{(A,A')}:j^*_g\cH_{A,A'}\isom\cH_{\sc_gA,\sc_gA'}
$$
whence a map
$\sd^g_{A,A'}:H_{A,A'}\xrightarrow{\cH_{A,A'}(g/X)}\cH_{A,A'}(g)
\xrightarrow{\tilde\mu^g_{(A,A'),\one_Y}}H_{\sc_gA,\sc_gA'}$ and
it is then easily seen that the rules : $A\mapsto\sc_gA$
for every $A\in\Ob(\cA_X)$ and $f\mapsto\sd^g_{A,A'}(f)$
for every morphism $f:A\to A'$ of $\cA^\sep_X$ define a functor
$$
\sd_g:\cA^\sep_X\to\cA^\sep_Y.
$$
Furthermore, a direct inspection yields a commutative diagram :
$$
\xymatrix@C+40pt{\Hom_{\cA_X}(A,A') \ar[r]^-{\tilde\omega_{A,A'}}
\ar[d]_{\sc_{g,(A,A')}} &
H_{A,A'} \ar[d]^{\sd^g_{A,A'}} \\
\Hom_{\cA_Y}(\sc_gA,\sc_gA') \ar[r]^-{\tilde\omega_{\sc_gA,\sc_gA'}} &
H_{\sc_gA,\sc_gA'}
}$$
where $\sc_{g,(A,A')}$ is the map given by the functor
$\sc_g:\cA_X\to\cA_Y$. In other words :
$$
\sd_g\circ\tilde\omega^X=\tilde\omega^Y\circ\sc_g.
$$
Next, let $Y'\xrightarrow{g'}Y\xrightarrow{g}X$ be two
morphisms of $\cC/X$; for every $A\in\Ob(\cA_X)$ there
exists a unique isomorphism of cartesian sections
\set\begin{equation}\label{eq_fart-2}
\beta^\blambda_{Y',\sc_{g'}\sc_gA}\isom\beta^\blambda_{Y',\sc_{gg'}A}
\qquad\text{such that}\qquad
\one_{Y'}\mapsto(\gamma^\sc_{(g,g'),A}:\sc_{g'}\sc_gA\isom\sc_{gg'}A)
\end{equation}
where $\gamma^\sc_{(\bullet,\bullet)}$ denotes the coherence
constraint of $\sc$. Denote by
$$
\gamma^\sd_{(g,g'),A}\in\Hom_{\cA_{Y'}^\sep}(\sc_g\sc_{g'}A,\sc_{gg'}A)
$$
the image of \eqref{eq_fart-2}, where the latter is seen as
an element of
$\cCart(\beta^\blambda_{Y',\sc_{g'}\sc_gA},\beta^\blambda_{Y',\sc_{gg'}A})(\one_{Y'})$.
A direct inspection of the construction yields a commutative diagram :
$$
\xymatrix@C+40pt{
(g\circ g')^\wedge_*\cCart(\beta^\blambda_{X,A},\beta^\blambda_{X,A'})
\ar[r]^-{g'^\wedge_*(\mu^g_{(A,A')})} \ar[d]_{\mu^{g\circ g'}_{(A,A')}} &
g'^\wedge_*\cCart(\beta^\blambda_{Y,\sc_gA},\beta^\blambda_{Y,\sc_gA'})
\ar[d]^{\mu^{g'}_{(\sc_gA,\sc_gA')}} \\
\cCart(\beta^\blambda_{Y',\sc_{gg'}A},\beta^\blambda_{Y',\sc_{gg'}A'})
\ar[r]^-{\delta^{g,g'}_{(A,A')}} &
\cCart(\beta^\blambda_{Y',\sc_{g'}\sc_gA},\beta^\blambda_{Y',\sc_{g'}\sc_gA'})
}$$
where $\delta^{g,g'}_{(A,A')}$ is the isomorphism of presheaves
\eqref{eq_pariah} induced by the isomorphisms \eqref{eq_fart-2}.
Let $(\delta^{g,g'}_{(A,A')})^a$ be the isomorphism of sheaves
on $C/X$ associated with $\delta^{g,g'}_{(A,A')}$; there follows a
commutative diagram :
$$
\xymatrix@C+20pt{
H_{A,A'} \ar[r]^-{\sd^g_{A,A'}} \ar[d]_{\sd^{gg'}_{A,A'}} &
H_{\sc_gA,\sc_gA'} \ar[d]^{\sd^{g'}_{\sc_gA,\sc_gA'}} \\
H_{\sc_{gg'}A,\sc_{gg'}A} \ar[r]^-{\sd^{g,g'}_{(A,A')}} &
H_{\sc_{g'}\sc_gA,\sc_{g'}\sc_gA'}
}$$
with $\sd^{g,g'}_{(A,A')}:=(\delta^{g,g'}_{(A,A')})^a_{\one_{Y'}}$.
This translates as the commutativity of the diagram :
$$
\qquad\qquad
{\spreaddiagramcolumns{+40pt}\diagram
\sc_{g'}\sc_gA \ar[r]^-{\sd_{g'}\circ\sd_g(f)}
\ar[d]_{\gamma^\sd_{(g,g'),A}} &
\sc_{g'}\sc_gA' \ar[d]^{\gamma^\sd_{(g,g'),A'}} \\
\sc_{gg'}A \ar[r]^-{\sd_{gg'}(f)} & \sc_{gg'}A'
\enddiagram}
\qquad
\text{for every $f\in\Hom_{\cA^\sep_X}(A,A')$}
$$
which means that the rule : $A\mapsto\gamma^\sd_{(g,g'),A}$
yields an isomorphism of functors
$\sd_{g'}\circ\sd_g\isom\sd_{gg'}$. Then it is easily seen
that the rules : $X\mapsto\cA^\sep_X$ and $g\mapsto\sd_g$
for every $X\in\Ob(\cC)$ and every morphism $g$ of $\cC$
yield a unital pseudo-functor whose coherence constraint
is the system of isomorphisms $\gamma^\sd_{\bullet\bullet}$ :
the detailed verification shall be left to the reader.
Likewise, it follows easily that the rule
$X\mapsto\tilde\omega^X$ defines a strict pseudo-natural
transformation $\tilde\omega^\bullet:\sc\Rightarrow\sd$. Set
$$
\cA^\sep:=\cFib(\sd)
$$
and let $\blambda^*$ be the distinguished cleavage of $\cA^\sep$,
whose associated pseudo-functor is naturally identified with
$\sd$ (see remark \ref{rem_distinguished-cleavage}(i)).
Then we have a unique cartesian functor
$$
\tilde\omega_\cA:\cA\to\cA^\sep
$$
restricting to $\tilde\omega^X:\cA_X\to\cA^\sep_X$ for every
$X\in\Ob(\cC)$. A direct inspection then shows that
$$
\beta^{\blambda^*}_{X,A}=\tilde\omega_\cA\circ\beta^\blambda_{X,A}
\qquad
\text{for every $X\in\Ob(\cC)$ and $A\in\Ob(\cA_X)$}
$$
and we have a natural isomorphism of presheaves :
$$
\cH_{A,A'}\isom\cCart(\beta^{\blambda^*}_{X,A},\beta^{\blambda^*}_{X,A'})
\qquad
\text{for every $A,A'\in\Ob(\cA_X)$}
$$
which identifies the morphism of presheaves induced by
$\tilde\omega_\cA$ as in \eqref{eq_F-pariah} :
$$
\cCart(\beta^\blambda_{X,A},\beta^\blambda_{X,A'})\to
\cCart(\tilde\omega_\cA\circ\beta^\blambda_{X,A},
\tilde\omega_\cA\circ\beta^\blambda_{X,A'})=
\cCart(\beta^{\blambda^*}_{X,A},\beta^{\blambda^*}_{X,A'})
$$
with the natural bicovering morphism of presheaves
$\cCart(\beta^\blambda_{X,A},\beta^\blambda_{X,A'})\to\cH_{A,A'}$.
Especially, $\cCart(\beta^{\blambda^*}_{X,A},\beta^{\blambda^*}_{X,A'})$
is a sheaf on $C/X$ for every such $X,A,A'$, so $\cA^\sep$ is
$1$-separated over $C$, by lemma \ref{lem_crit-0-1_separation}.
Consider now any $1$-separated fibration $\cA'$ over $C$,
and any cartesian functor $F:\cA\to\cA'$. For every
$X\in\Ob(\cC)$ and $A,A'\in\cA_X$, the morphism of presheaves
\set\begin{equation}\label{eq_frufru}
\cCart(\beta^\blambda_A,\beta^\blambda_{A'})\to
\cCart(F\circ\beta^\blambda_A,F\circ\beta^\blambda_{A'})
\end{equation}
as in \eqref{eq_F-pariah} factors uniquely through $\cH_{A,A'}$,
since $\cCart(F\circ\beta^\blambda_A,F\circ\beta^\blambda_{A'})$
is a sheaf on $C/X$ (lemma \ref{lem_crit-0-1_separation}).
There follows a map $F^\sep_{A,A'}:H_{A,A'}\to\Hom_{\cA'}(FA,FA')$,
whose composition with $\tilde\omega_{A,A'}$ equals
$F_{A,A'}:\Hom_\cA(A,A')\to\Hom_{\cA'}(FA,FA')$. It is easily
seen that the rules : $A\mapsto F^\sep A:=FA$ and
$(A,A')\mapsto F^\sep_{A,A'}$ yield a well defined functor
$F^\sep_{|X}:\cA^\sep_X\to\cA'_X$ such that
$F^\sep_{|X}\circ\tilde\omega^X$ agrees with the restriction
$\cA_X\to\cA'_X$ of $F$.

Let $\blambda'$ be a unital cleavage for $\cA'$, and $\sc'$
its associated pseudo-functor; then $F$ corresponds to
a pseudo-natural transformation $\alpha:\sc\Rightarrow\sc'$,
and denote by $\tau^\alpha_\bullet$ the coherence constraint
of $\alpha$. For every $X\in\Ob(\cC)$ and $A\in\Ob(\cA_X)$
we have also a unique isomorphism
\set\begin{equation}\label{eq_worry}
\beta^{\blambda'}_{FA}\isom F\circ\beta^\blambda_A
\qquad\text{such that}\qquad
\one_X\mapsto\one_{FA}.
\end{equation}
Explicitly, \eqref{eq_worry} assigns to every
$(g:Y\to X)\in\Ob(\cC/X)$ the isomorphism
$\tau^\alpha_{g,A}:\sc'_gFA\isom F\sc_gA$ of $\cA'_Y$.
On the other hand, we have a morphism of presheaves
\set\begin{equation}\label{eq_please}
\cCart(\beta^\blambda_A,\beta^\blambda_A)\to
\cCart(\beta^{\blambda'}_{FA},\beta^{\blambda'}_{FA})
\end{equation}
described explicitly as follows. For every
$(g:Y\to X)\in\Ob(\cC/X)$, every morphism of cartesian
sections $t_\bullet:\beta^\blambda_A\circ g_*\Rightarrow
\beta^\blambda_{A'}\circ g_*$ is determined by $t_{\one_Y}$ :
namely, $t_h=\sc_h(t_{\one_Y})$ for every $h\in\Ob(\cC/Y)$;
then \eqref{eq_please} maps every such $t_\bullet$ to
the morphism
$t'_\bullet:\beta^{\blambda'}_{FA}\circ g_*\Rightarrow
\beta^{\blambda'}_{FA'}\circ g_*$ with $t'_h:=\sc'_h(Ft)$
for every $h\in\Ob(\cC/Y)$. With this notation, the
condition that $\tau^\alpha_g$ is a natural isomorphism
$\sc'_g\circ F\isom F\circ\sc_g$ for every $g\in\Ob(\cC/X)$
translates as the commutativity of the diagram of presheaves
on $\cC/X$ :
$$
\xymatrix{ \cCart(\beta^\blambda_A,\beta^\blambda_A)
\ar[r] \ar[d] &
\cCart(\beta^{\lambda'}_{FA},\beta^{\blambda'}_{FA'})
\ar[d] \\
\cCart(F\circ\beta^\blambda_A,F\circ\beta^\blambda_A) \ar[r] &
\cCart(F\circ\beta^\blambda_A,\beta^{\lambda'}_{FA})
}$$
whose top horizontal (resp. left vertical) arrow is
\eqref{eq_please} (resp. \eqref{eq_frufru}) and whose
other two arrows are the morphisms \eqref{eq_pariah}
induced by the isomorphisms \eqref{eq_worry} (and by
the identities of $\beta^{\lambda'}_{FA}$ and of
$F\circ\beta^\blambda_A$). After taking associated sheaves,
we obtain a similar diagram, in which
$\cCart(\beta^\blambda_A,\beta^\blambda_A)$ is replaced by
$\cH_{A,A'}$.
Unwinding the definitions, we find that the commutativity
of the latter diagram yields the commutativity of the following
diagram of morphisms of $\cA_Y$ :
$$
{\spreaddiagramcolumns{+20pt}\diagram
\sc'_gFA \ar[r]^-{\sc'_gF^\sep f} \ar[d]_{\tau^\alpha_{g,A}} &
\sc'_gFA' \ar[d]^{\tau^\alpha_{g,A'}} \\
F\sc_gA \ar[r]^-{F^\sep\sd_gf} & F\sc_gA'
\enddiagram}
\qquad\qquad
\text{for every $f\in H_{A,A'}$}.
$$
Thus, $\tau^\alpha_g$ extends to a natural isomorphism
$F^\sep_{|Y}\circ\sc_g\isom\sc'_g\circ F^\sep_{|X}$, and
then it is clear that the rule $X\mapsto F^\sep_{|X}$ yields
a pseudo-natural transformation $\alpha':\sd\Rightarrow\sc'$
whose coherence constraint is again $\tau^\alpha_\bullet$.
Finally, $\alpha'$ corresponds to a unique cartesian
functor $F^\sep:\cA^\sep\to\cA'$ whose restriction
$\cA^\sep_X\to\cA'_X$ agrees with $F^\sep_{|X}$ for every
$X\in\Ob(\cC)$. It is then clear that $F^\sep$ is the
unique such cartesian functor whose composition with
$\tilde\omega_\cA$ equals $F$. This universal property
then easily implies that the rules : $\cA\mapsto\cA^\sep$
and $F\mapsto F^\sep$ yield a left adjoint $(-)^\sep$ for
the inclusion functor $1\tdu\PreStack(C)\Fib(\cC)$ : the
details shall be left to the reader.

Lastly, let $\mu:F\Rightarrow G$ be any natural
$\cC$-transformation of cartesian functors $F,G:\cA\to\cA'$
between $\cC$-fibrations. Then for every $X\in\Ob(\cC)$
and every $A\in\cA_X$, the morphism $\mu_A:FA\to GA$
determines a unique natural transformation
$\mu^*_A:F\circ\beta^\blambda_A\Rightarrow G\circ\beta^\blambda_A$
such that $\one_X\mapsto\mu_A$. There follows for every
$A,A'\in\Ob(\cA)$ a commutative diagram of presheaves :
$$
\xymatrix{ \cCart(\beta^\blambda_A,\beta^\blambda_{A'})
\ar[r] \ar[d] &
\cCart(F\circ\beta^\blambda_A,F\beta^\blambda_{A'}) \ar[d] \\
\cCart(G\circ\beta^\blambda_A,G\circ\beta^\blambda_{A'})
\ar[r] & \cCart(F\circ\beta^\blambda_A,G\circ\beta^\blambda_{A'})
}$$
whose top horizontal and left vertical arrows are
the morphisms \eqref{eq_frufru}, and whose other two
arrows are the morphisms \eqref{eq_pariah} induced
by $\mu^*_A$ and $\mu^*_{A'}$. After taking associated
sheaves on $C/X$, we deduce a commutative diagram
$$
\xymatrix{ \Hom_{\cA^\sep}(A,A') \ar[r] \ar[d] &
\Hom_{\cB^\sep}(FA,FB) \ar[d] \\
\Hom_{\cA^\sep}(GA,GA') \ar[r] & \Hom_{\cA^\sep}(FA,GA')  
}$$
which shows that $\mu$ is also a natural transformation
$F^\sep\Rightarrow G^\sep$. Clearly the rule $\mu\mapsto\mu^\sep$
is inverse to the rule $(\nu:F^\sep\Rightarrow G^\sep)\mapsto
\nu*\tilde\omega_\cA$, so $(-)^\sep$ is the sought strong
left $2$-adjoint of the inclusion pseudo-functor.
\end{proof}

\sset\subsubsection{}\label{subsec_+construction}
Let $(\cC,J)$ be a small site. Our next task is to construct
the analogue for stacks of the functor $F\mapsto F^+$ on presheaves
(see \eqref{subsec_asso-topoi}). To this aim, recall
that for every $X\in\Ob(\cC)$ the set $J(X)$ of sieves covering $X$
is cofiltered, and if $\cS'\subset\cS$ are two such sieves and
$F:\cE\to\cC$ any fibration, the inclusion functor $i:\cS'\to\cS$
induces a functor
$$
\sCart_\cC(i,\cE):\sCart_\cC(\cS,\cE)\to\sCart_\cC(\cS',\cE)
$$
whence a well defined filtered system of small categories
$C_X:=(\sCart_\cC(\cS,\cE)~|~\cS\in J(X)^o)$ associated with
every object $X$ of $\cC$. If $\cE$ has small fibres, we may
then consider the functor
$$
\sc_F^+:\cC^o\to\bCat
\qquad
X\mapsto\colim_{\cS\in J(X)^o}\sCart_\cC(\cS,\cE)
$$
where the colimit is explicitly given by the construction
detailed in example \ref{ex_fil-colim-in-Cat}. For any
morphism $f:Y\to X$ in $\cC$, the corresponding functor
$\sc_F^+(f):\sc_F^+(X)\to\sc_F^+(Y)$ is obtained as follows.
For every $\cS\in J(X)$, the functor
$f_{*|\cS}:\cS\times_Xf\to\cS$ induces a functor
$$
\sCart_\cC(f_{*|\cS},\cE):
\sCart_\cC(\cS,\cE)\to\sCart(\cS\times_Xf,\cE)\to\sc_F^+(Y)
$$
(see remark \ref{rem_sieves-and-sub}(v)) and the system
of such functors obviously forms a cocone of vertex $\sc_F^+(Y)$
and basis $C_X:J(X)^o\to\bCat$. There follows a unique functor
$\sc_F^+(f):\sc_F^+(X)\to\sc_F^+(Y)$ whose composition with the
natural functor $\sCart_\cC(\cS,\cE)\to\sc_F^+(X)$ agrees with
$\sCart_\cC(f_{*|\cS},\cE)$ for every $\cS\in J(X)$. It is then
easily seen that
$$
\sc_F^+(f\circ g)=\sc_F^+(g)\circ\sc_F^+(f)
\qquad
\text{for every pair of morphisms
$Z\xrightarrow{\ g\ }Y\xrightarrow{\ f\ }X$ in $\cC$}.
$$
With the notation of \eqref{subsec_fib-from-pseudo}, we set
$$
\cE^+:=\cFib(\sc^+_F)
$$
and let $F^+:\cE^+\to\cC$ be the structure functor of $\cE^+$.
Recall also that $\cE^+$ is endowed with a distinguished
split cleavage $\blambda^+_F$ whose associated pseudo-functor
is naturally identified with $\sc^+_F$ (see remark
\ref{rem_distinguished-cleavage}(i)). Thus, the objects of
$\cE^+$ are represented by the pairs
$$
[X,\psi:\cS\to\cE]
$$
with $X\in\Ob(\cC)$ and $\psi$ is a $\cC$-cartesian functor
defined on some covering sieve $\cS\in J(X)$. A morphism
$[f,\sigma]:[X',\psi':\cS'\to\cE]\to[X,\psi:\cS\to\cE]$ is
represented by the datum of a morphism $f:X'\to X$ in $\cC$,
and a natural $\cC$-transformation
$\sigma:\psi'_{|\cS''}\Rightarrow(\psi\circ f_*)_{|\cS''}$, where
$\cS''\subset\cS'\cap(\cS\times_Xf)$ is a sieve covering $X'$.

\sset\subsubsection{}\label{subsec_functoriality-of-+}
Next, let $F':\cE'\to\cC$ be another fibration with small fibres,
and $G:\cE\to\cE'$ any $\cC$-cartesian functor; the system of
induced functors
$$
\sCart_\cC(\cS,G):\sCart_\cC(\cS,\cE)\to\sCart_\cC(\cS,\cE')
\qquad
\text{for every $X\in\Ob(\cC)$ and $\cS\in J(X)$}
$$
yields, after taking colimits, a strict pseudo-natural
transformation
$$
\sc^+_G:\sc^+_F\Rightarrow\sc^+_{F'}
\qquad\text{and then a $\cC$-cartesian functor :}
\qquad
G^+:=\cFib(\sc^+_G):\cE^+\to\cE'^+.
$$
Moreover, if $\omega:G\Rightarrow G'$ is any natural
$\cC$-transformation of $\cC$-cartesian functors
$G,G':\cE\to\cE'$, we get a system of natural transformations
$$
\sCart_\cC(\cS,\omega):\sCart_\cC(\cS,G)\Rightarrow\sCart_\cC(\cS,G')
\qquad
\text{for every $X\in\Ob(\cC)$ and $\cS\in J(X)$}
$$
whence again, after taking colimits, a modification
$$
\sc_\omega^+:\sc^+_G\leadsto\sc^+_{G'}
\qquad
\text{whence a natural transformation :}
\qquad
\omega^+:=\Fib(\sc^+_\omega):G^+\Rightarrow G'^+.
$$
It is then clear that the rules :
$(F:\cE\to\cC)\mapsto F^+$, $(G:\cE\to\cE')\mapsto G^+$ and
$(\omega:G\Rightarrow G')\mapsto\omega^+$ yield a well
defined strict pseudo-functor
$$
(-)^+:\Fib(\cC)\to\Fib(\cC).
$$
Moreover, since $\cC/X\in J(X)$ for every $X\in\Ob(\cC)$,
we have an obvious strict pseudo-natural transformation of
strict pseudo-functors :
$$
\sj_\cE:\sCart_\cC(\cC/-,\cE)\Rightarrow\sc^+_\cE
$$
(notation of \eqref{subsec_we-mention}), whence a $\cC$-cartesian
functor
$$
j_\cE:=\Fib(\sj_F):\sC(\cE)\to\cE^+
$$
and it is clear that the rule $(F:\cE\to\cC)\mapsto j_\cE$ yields
a strict pseudo-natural transformation of strict pseudo-functors.

\begin{remark}\label{rem_cartesian-to-E+}
Let $F:\cE\to\cC$ be any fibration with small fibres,
$X\in\Ob(\cC)$ any object, and $\cS$ any sieve of $\cC/X$
covering $X$. By unwinding the definitions, and taking into
account lemma \ref{lem_distinguished-cleavage}(i), we see that
a $\cC$-cartesian functor $\cS\to\cE^+$ is the datum of :
\begin{itemize}
\item
for every object $Y\xrightarrow{f}X$ of $\cS$, a sieve
$\cS^{(f)}$ covering $Y$ and a $\cC$-cartesian functor
$$
\sigma^{(f)}:\cS^{(f)}\to\cE
$$
\item
for every morphism
$h/X:(Y'\xrightarrow{f'}X)\to(Y\xrightarrow{f}X)$ of $\cS$,
a sieve $\cS^{(h/X)}\subset\cS^{(f')}\cap(\cS^{(f)}\times_Yh)$
covering $Y'$ and a natural isomorphism of $\cC$-functors
$$
\sigma^{(h/X)}:\sigma^{(f')}_{|\cS^{(h/X)}}\isom
(\sc^+_F(h)(\sigma^{(f)}))_{|\cS^{(h/X)}}
=(\sigma^{(f)}\circ h_*)_{|\cS^{(h/X)}}
$$
where $h_*:\cS^{(f)}\times_Yh\to\cS^{(f)}$ is the functor
$(Z\xrightarrow{g}Y')\mapsto(Z\xrightarrow{h\circ g}Y)$
\item
for every $h/X$ as in the foregoing and every other morphism
$h'/X:(Y''\xrightarrow{f''}X)\to(Y'\xrightarrow{f}X)$ of
$\cS_\sigma$, a covering sieve
$\cS^{(h,h'/X)}\subset\cS^{(h'/X)}\cap(\cS^{(h/X)}\times_{Y'}h')$
such that
\set\begin{equation}\label{eq_one-last-effort}
(\sigma^{(h/X)}*h'_*)_{|\cS^{(h,h'/X)}}\odot\sigma^{(h'/X)}_{|\cS^{(h,h'/X)}}=
\sigma^{(h\circ h'/X)}_{|\cS^{(h,h'/X)}}.
\end{equation}
\end{itemize}
\end{remark}

\begin{lemma}\label{lem_pair-of-sieves}
Let $i\leq 2$ be an integer, $F:\cE\to\cC$ an $i$-separated
prestack, $X\in\Ob(\cC)$. Then every inclusion $\cT\subset\cS$
of sieves of\/ $\cC/X$ covering $X$ induces an $i$-faithful functor
$$
j:\sCart_\cC(\cS,\cE)\to\sCart_\cC(\cT,\cE).
$$
\end{lemma}
\begin{proof} It is an immediate consequence of lemma
\ref{lem_was-pair-of-sieves}(i).
\end{proof}

\begin{lemma}\label{lem_bootstrap}
Let $\cC$ be a small category, $F:\cE\to\cC$ a fibration
with small fibres. We have :
\begin{enumerate}
\item
The prestack $F^+:\cE^+\to\cC$ is $0$-separated.
\item
If $F$ is a $0$-separated prestack, then $F^+$ is a
$1$-separated prestack.
\item
If $F$ is a $1$-separated prestack, then $F^+$ is a stack.
\item
If $F$ is a stack, the functor $j_\cE:\sC(\cE)\to\cE^+$
is an equivalence.
\item
There exists an essentially commutative diagram of strict
pseudo-functors (notation of definition {\em\ref{def_sheaf}(ii)}) :
$$
\xymatrix{ \cC^\wedge \ar[rr]^-\cFib \ar[d]_{(-)^+} & &
0\tdu\PreStack(\cC,J) \ar[d]^{(-)^+} \\
(\cC,J)^\sep \ar[rr]^-\cFib & & 1\tdu\PreStack(\cC,J).
}$$
\end{enumerate}
\end{lemma}
\begin{proof} Let $X\in\Ob(\cC)$ be any object, and
$\cT$ any covering sieve of $\cC/X$; let also
$[\phi],[\psi]\in(F^+)^{-1}X$, and $\cS$ a sieve covering
$X$ such that $[\phi]$ and $[\psi]$ are represented by
$\cC$-cartesian functors $\phi,\psi:\cS\to\cE$.

(i): Let $[\alpha],[\beta]:[\phi]\to[\psi]$ be morphisms
of $(F^+)^{-1}X$ whose images agree under the functor
$$
j:(F^+)^{-1}X\to\sCart_\cC(\cT,\cE^+).
$$
deduced from the split cleavage $\sc^+_F$.
We need to check that $\alpha=\beta$. To this aim, may
assume that $[\alpha]$ and $[\beta]$ are represented by
natural $\cC$-transformations
$\alpha,\beta:\phi\Rightarrow\psi$; then, for every
$(Y\xrightarrow{f}X)\in\Ob(\cT)$ there exists a sieve
$\cT^{(f)}\subset\cS\times_Xf$ covering $Y$, such that
$$
(\alpha*f_*)_{|\cT^{(f)}}=(\beta*f_*)_{|\cT^{(f)}}.
$$
Now, let $\cS'\subset\cS$ be the sieve generated by
$\bigcup_{f\in\Ob(\cT)}f_*(\Ob(\cT^{(f)}))$; by remark
\ref{rem_topology}(iii), the sieve $\cS'$ covers
$X$, and clearly $\alpha_{|\cS'}=\beta_{|\cS'}$, whence
the contention.

(ii): Let $\alpha:j([\phi])\Rightarrow j([\psi])$ be
a given natural $\cC$-transformation; we need to check
that $\alpha=j(\alpha')$ for some sieve $\cS'\subset\cS$
covering $X$ and some natural $\cC$-transformation
$\alpha':\phi_{|\cS'}\Rightarrow\psi_{|\cS'}$. However,
$\alpha$ is described by a datum as follows :
\begin{itemize}
\item
For every $(Y\xrightarrow{f}X)\in\Ob(\cT)$, a sieve
$\cT^{(f)}\subset\cS\times_Xf$ covering $Y$ and a
natural $\cC$-transformation
$\alpha^{(f)}:(\phi\circ f_*)_{|\cT^{(f)}}\Rightarrow
(\psi\circ f_*)_{|\cT^{(f)}}$
\item
for every morphism
$(h/X):(Y'\xrightarrow{f'}X)\to(Y\xrightarrow{f}X)$ in
$\cT$, a sieve $\cT^{(h/X)}$ covering $Y'$ and contained
in $\cT^{(f')}\cap(\cT^{(f)}\times_Yh)$ such that
$$
\alpha^{(f')}_{|\cT^{(h/X)}}=(\alpha^{(f)}*h_*)_{|\cT^{(h/X)}}.
$$
\end{itemize}
But then, in view of lemma \ref{lem_pair-of-sieves}, we may
take already $\cT^{(h/X)}:=\cT^{(f')}\cap(\cT^{(f)}\times_Yh)$.
With this notation, let $\cS'$ be the sieve generated by
$\bigcup_{f\in\Ob(\cT)}f_*(\Ob(\cT^{(f)}))$; then $\cS'$ covers
$X$, by remark \ref{rem_topology}(iii). For every
$(Y\xrightarrow{g}X)\in\Ob(\cS')$, pick any $f\in\Ob(\cT)$
such that $g=g'\circ f$ for some $g'\in\Ob(\cT^{(f)})$, and set
$$
\alpha'_g:=\alpha^{(f)}_{g'}:\phi(g)\to\psi(g).
$$
Let us check that this definition is independent of the
choices of $f$ and $g'$. Indeed, suppose that $g=g''\circ f'$
for some other $f'\in\Ob(\cT)$ and some $g''\in\Ob(\cT^{(f')})$;
then we have
$$
(\alpha^{(f)}*g'_*)_{|\cT'}=(\alpha^{(g)})_{|\cT'}=
(\alpha^{(f')}*g''_*)_{|\cT'}
\quad
\text{with
$\cT':=\cT^{(g)}\cap(\cT^{(f)}\times_Yg')\cap(\cT^{(f')}\times_Yg'')$}
$$
and by invoking again lemma \ref{lem_pair-of-sieves}, we
deduce that
$$
(\alpha^{(f)}*g'_*)_{|\cT''}=(\alpha^{(f')}*g''_*)_{|\cT''}
\qquad
\text{with $\cT'':=(\cT^{(f)}\times_Yg')\cap(\cT^{(f')}\times_Yg'')$}.
$$
Especially, we get : $\alpha^{(f)}_{g'}=(\alpha^{(f)}*g'_*)_{\one_Y}=
(\alpha^{(f')}*g''_*)_{\one_Y}=\alpha^{(f')}_{g''}$. It is then clear
that $\alpha'$ is a well defined natural $\cC$-transformation
as sought.

(v) follows now from (ii) and a direct inspection of the
construction of the functors $(-)^+$ for presheaves and for
prestacks.

(iii): Let $\sigma:\cT\to\cE^+$ be a $\cC$-cartesian functor;
hence, this is a datum
$$
\Bigl(\!(\cS^{(f)},\sigma^{(f)}~|~f\!\in\!\Ob(\cT)),
(\cS^{(h/X)},\sigma^{(h/X)}~|~\text{$f'\xrightarrow{h}f$
in $\cT$}),
(\cS^{(h,h'/X)}~|~\text{$f''\xrightarrow{h'}f'\xrightarrow{h}f$
in $\cT$})\!\Bigr)
$$
as in remark \ref{rem_cartesian-to-E+}. But then, lemma
\ref{lem_pair-of-sieves} says that we may even take :
$$
\cS^{(h/X)}:=\cS^{(f')}\cap(\cS^{(f)}\times_Yh)
\qquad
\cS^{(h,h'/X)}:=\cS^{(f'')}\cap(\cS^{(f')}\times_{Y'}h')
\cap(\cS^{(f)}\times_Y(h\circ h'))
$$
for every $f\in\Ob(\cT)$, every morphism $f'\xrightarrow{h}f$
in $\cT$, and every pair $f''\xrightarrow{h'}f'\xrightarrow{h}f$
of morphisms in $\cT$. Let $\cT'\subset\cC/X$ be the sieve
generated by $\bigcup_{f\in\Ob(\cT)}f_*(\Ob(\cS^{(f)}))$, which
covers $X$, by remark \ref{rem_topology}(iii). We consider
the functor
$$
G_{\cT'}:\cT'\to\bCat/(\cC/X)
\qquad
(Y\xrightarrow{g}X)\mapsto(\cC/X)/g=\cC/Y
$$
defined as in \eqref{subsec_towards-top}, and we construct
a pseudo-cocone
$$
\psi_\bullet:G_{\cT'}\Rightarrow\sF_\cE
$$
as follows. For every $(Y\xrightarrow{g}X)\in\Ob(\cT')$,
choose $g'\in\Ob(\cT)$ and $g''\in\Ob(\cS^{(g')})$ with
$$
g=g'\circ g''
\qquad\text{and set :}\qquad
\psi_g:=\sigma^{(g')}\circ g''_*:\cC/Y\to\cE.
$$
Next, let $g_1,g_2\in\Ob(\cT')$ and
$h/X:(Z_1\xrightarrow{g_1}X)\to(Z_2\xrightarrow{g_2}X)$ a
morphism of $\cT'$; with this notation, we have the
$\cC$-cartesian functors
$$
\sigma^{(g'_1)}\circ g''_{1*},\sigma^{(g'_2)}\circ(g''_2\circ h)_*:
\cC/Z_1\to\cE
$$
and two isomorphisms of $\cC$-functors
$$
(\sigma^{(g'_1)}\circ g''_{1*})_{|\cS^{(g_1)}}
\xLeftarrow{\ \sigma^{(g''_1/X)}\ }
\sigma^{(g_1)}\xRightarrow{\ \sigma^{(g''_2\circ h/X)}\ }
(\sigma^{(g'_2)}\circ(g''_2\circ h)_*)_{|\cS^{(g_1)}}.
$$
Since $F$ is $1$-separated, it follows that the
composition $\sigma^{(g''_2\circ h/X)}_{|\cS^{(g_1)}}\odot
(\sigma^{(g''_1/X)})_{|\cS^{(g_1)}}^{-1}$ extends uniquely
to an isomorphism of $\cC$-functors
$$
\tau^\psi_{h/X}:\psi_{g_1}\isom\psi_{g_2}\circ h_*.
$$
We claim that the rule $g\mapsto\psi_g$ yields a pseudo-cocone
as sought, with coherence constraint given by the system of
isomorphisms $\tau^\psi_\bullet$. Indeed, let
$h'/X:(Z_0\xrightarrow{g_0}X)\to(Z_1\xrightarrow{g_1}X)$ be
another morphism of $\cT'$; then we have as well the isomorphisms
of $\cC$-functors
$$
\xymatrix{ (\sigma^{(g'_0)}\circ g''_{0*})_{|\cS^{(g_0)}}
& \ar@{=>}[l]_-{\sigma^{(g''_0/X)}} \sigma^{(g_0)} 
\ar@{=>}[r]^-{\sigma^{(g''_1\circ h'/X)}}
\ar@{=>}[d]^-{\sigma^{(g''_2\circ h'\circ h/X)}} &
(\sigma^{(g'_1)}\circ(g''_1\circ h')_*)_{|\cS^{(g_0)}} \\
& (\sigma^{(g'_2)}\circ(g''_2\circ h\circ h')_*)_{|\cS^{(g_0)}}
}$$
and with this notation, we have :
$$
(\tau^\psi_{h'/X})_{|\cS^{(g_0)}}=
\sigma^{(g''_1\circ h'/X)}\odot(\sigma^{(g''_0/X)})^{-1}
\qquad
(\tau^\psi_{h\circ h'/X})_{|\cS^{(g_0)}}=
(\sigma^{(g''_2\circ h\circ h'/X)})\odot(\sigma^{(g''_0/X)})^{-1}.
$$
Set $\cS:=\cS^{(g''_2\circ h,h'/X)}\cap\cS^{(g''_1,h'/X)}\subset\cS^{(g_0)}$.
We claim that
$$
(\tau^\psi_{h\circ h'/X})_{|\cS}=
(\tau^\psi_{h/X}*h'_*)_{|\cS}\odot
(\tau^\psi_{h'/X})_{|\cS}.
$$
Indeed, the latter identity follows by combining the
following ones :
$$
\begin{aligned}
\sigma^{(g''_2\circ h\circ h'/X)}_{|\cS}=&\,
(\sigma^{(g''_2\circ h/X)}*h'_*)_{|\cS}\odot\sigma^{(h'/X)}_{|\cS} \\
\sigma^{(g''_1\circ h'/X)}_{|\cS}=&\,
(\sigma^{(g''_1/X)}*h'_*)_{|\cS}\odot\sigma^{(h'/X)}_{|\cS}
\end{aligned}
$$
that are provided by remark \ref{rem_cartesian-to-E+}.
Since $\cS$ covers $Z_0$ and $F$ is $1$-separated,
it follows that
$$
\tau^\psi_{h\circ h'/X}=(\tau^\psi_{h/X}*h'_*)\odot\tau^\psi_{h'/X}
$$
which shows that $\tau^\psi$ satisfies the required
coherence axioms. By lemma \ref{lem_new-pseudo-col}, we
deduce that there exist a $\cC$-functor $\psi:\cT'\to\cE$
and an invertible modification
$\Xi:\sF_\psi\odot\hat\eps_{\cT'}\leadsto\psi_\bullet$, where
$\hat\eps_{\cT'}:G_{\cT'}\Rightarrow\sF_{\cT'}$ is the universal
pseudo-cocone defined as in \eqref{subsec_towards-top};
explicitly, this amounts to a system of oriented squares
with invertible $2$-cells :
$$
{\diagram \cC/Y \ar[d]_{g''_*} \ar[r]^-{g_*}
\drtwocell\omit{_\ \ \ \ \Xi_g} &
\cT' \ar[d]^\psi \\
\cS^{(g')} \ar[r]_-{\sigma^{(g')}} & \cE
\enddiagram}
\qquad
\text{for every $g\in\Ob(\cT')$}
$$
from which it follows easily that $\psi$ is $\cC$-cartesian
(here $g=g'\circ g''$ is the factorisation used to define
$\psi_g$). The compatibility condition for $\Xi$ amounts
to the identity :
\set\begin{equation}\label{eq_compat-for-Xi}
\Xi_{g_2}*h_*=\tau^\psi_{h/X}\odot\Xi_{g_1}
\qquad
\text{for every morphism
$(Z_1\xrightarrow{g_1}X)\xrightarrow{h/X}(Z_2\xrightarrow{g_2}X)$
of $\cT'$}.
\end{equation}
Thus, the functor $\psi$ represents an object $[\psi]$ of the
fibre category $\cE^+(X)$, and to conclude, we need to check
that $j([\psi])$ is isomorphic to $\sigma$. However, since
we know already from (ii) that $\cE^+$ is a $1$-separated
prestack, by lemma \ref{lem_pair-of-sieves} it suffices to
find an isomorphism between the images of $j([\psi])$ and
$\sigma$ in $\sCart_\cC(\cT',\cE^+)$. This comes down to
exhibiting isomorphisms of functors
$$
\theta_g:\psi\circ g_*\isom\sigma^{(g)}
\qquad
\text{for every $g\in\Ob(\cT')$}
$$
(where $g_*:\cS^{(f)}\to\cT'$ is the usual functor with
$t\mapsto g\circ t$ for every $t\in\Ob(\cS^{(f)})$) such that
\set\begin{equation}\label{eq_last-effort}
\sigma^{(h/X)}\odot\theta_{g_1|\cS^{(h/X)}}=
(\theta_{g_2}*h_*)_{|\cS^{(h/X)}}
\qquad
\text{for every morphism $g_1\xrightarrow{h/X}g_2$ of $\cT'$}.
\end{equation}
To this aim, we set :
$$
\theta_g:=(\sigma^{(g''/X)})^{-1}\odot(\Xi_g)_{|\cS^{(g)}}.
$$
Since $\cE$ is $1$-separated, the identities \eqref{eq_last-effort}
can be checked after restriction to any covering sieve contained
in $\cS^{(h/X)}$; but on such a suitable sieve we may apply the
identities \eqref{eq_one-last-effort}, and combining with
\eqref{eq_compat-for-Xi} we conclude easily : the details
shall be left to the reader.

(iv): It suffices to check that $j_\cE$ is a fibrewise
equivalence when $\cE$ is a stack (corollary
\ref{cor_fibrations}(i)), and this follows by direct
inspection of the definitions.
\end{proof}

\begin{theorem}\label{th_stackfication}
For every small site $(\cC,J)$, the inclusion strict
pseudo-functor
$$
\sF:\Stack(\cC,J)\to\Fib(\cC)
$$
admits a left $2$-adjoint
$$
\Fib(\cC)\to\Stack(\cC,J)
\qquad
(F:\cE\to\cC)\mapsto(F^a:\cE^a\to\cC)
$$
which assigns to every fibration over $\cC$ its
{\em associated stack}. Moreover, we have a
pseudo-commutative diagram of pseudo-functors
(see remark {\em\ref{rem_pseudo-commutative}}) :
\set\begin{equation}\label{eq_free-palestine}
{\diagram \cC^\wedge \ar[rr]^-\cFib \ar[d]_{(-)^a} & &
\Fib(\cC) \ar[d]^{(-)^a} \\
(\cC,J)^\sim \ar[rr]^-\cFib & & \Stack(\cC,J)
\enddiagram}
\end{equation}
whose left vertical arrow is the usual left adjoint to
the inclusion functor for presheaves.
\end{theorem}
\begin{proof} Recall first that the counit of the $2$-adjoint
pair $(\sF,\sC)$ of \eqref{subsec_we-mention} is a
$\cC$-equivalence  $\sev^\cE:\sC(\cE)\to\cE$ for every
fibration $F:\cE\to\cC$ (theorem \ref{th_split-fibration}).
Moreover, the unit of the same $2$-adjoint pair assigns to
every split fibration $(F,\blambda)$ a $\cC$-equivalence
$\eta_\cE:\cE\to\sC(\cE)$, deduced from the strict
pseudo-natural equivalence $\beta^\blambda$ defined in the
proof of theorem \ref{th_split-fibration} (see remark
\ref{rem_mathsfev}(i)). There follows a $\cC$-equivalence :
$$
\omega_\cE:=:\sC(\cE)^+\xrightarrow{\ (\sev^\cE)^+\ }\cE^+
\xrightarrow{\ \eta_{\cE^+}\ }\sC(\cE^+).
$$
By inspecting the definitions, we see that $\omega_\cE$
assigns to every object $[X,\phi:\cS\to\sC(\cE)]$ of $\sC(\cE)^+$
the cartesian section $\omega_\cE([X,\phi])\in\cE^+(X)$
given by the rule :
$$
(f:Y\to X)\mapsto(\cS\times_Xf\xrightarrow{\ f_*\ }\cS
\xrightarrow{\ \sev^\cE\circ\phi\ }\cE)
$$
and notice that
$\sev^\cE\circ\phi\circ f_*(h)=\phi_{f\circ h}(\one_Z)$ for every
$(h:Z\to Y)\in\Ob(\cS\times_Xf)$. For every morphism
$g/X:(f':Y'\to X)\to(f:Y\to X)$ in $\cC/X$, the induced
natural transformation
$\omega_\cE([X,\phi])(g):\omega_\cE([X,\phi])(f')\Rightarrow
\omega_\cE([X,\phi])(f)$ is the identity
$\sev^\cE\circ\phi\circ f'_*\isom\sev^\cE\circ\phi\circ f_*\circ g_*$.

\begin{claim}\label{cl_again-plus-commutes}
(i)\ \
For every two $\cC$-fibrations
$\cA\xrightarrow{\phi}\cC\xleftarrow{\phi'}\cA'$, the functor
$$
\sC_{\cA,\cA'}:\sCart_\cC(\cA,\cA')\to\sCart_\cC(\sC(\cA),\sC(\cA'))
\qquad
F\mapsto\sC(F)
\qquad
(\alpha:F\Rightarrow F')\mapsto
\sC(\alpha)
$$
is an equivalence.

(ii)\ \
For every $\cC$-fibration $\cA\to\cC$ the following diagram
of functors commutes:
$$
\xymatrix{ \sCart_\cC(\sC(\cE),\cA) \ar[rrr]^-{(-)^+}
\ar[d]_{\sC_{\sC(\cE),\cA}} & & &
\sCart_\cC(\sC(\cE)^+,\cA^+) \ar[d]^{\sCart_\cC(j_{\sC(\cE)},\cA^+)} \\
\sCart_\cC(\sC(\sC(\cE)),\sC(\cA))
\ar[rrr]^-{\sCart_\cC(\sC(\sC(\cE)),j_\cA)} & & &
\sCart_\cC(\sC(\sC(\cE)),\cA^+).
}$$
\end{claim}
\begin{pfclaim} Assertion (ii) follows by direct inspection.
For (i) let us notice the commutative diagram :
$$
\xymatrix{ \sCart_\cC(\cA,\cA') \ar[drr]_{\sCart_\cC(\sev^\cA,\cA')\ \ \ }
\ar[rrrr]^-{\sC_{\cA,\cA'}} & & & &
\sCart_\cC(\sC(\cA),\sC(\cA')) \ar[dll]^{\sCart_\cC(\cA,\sev^{\cA'})} \\
& & \sCart_\cC(\sC(\cA),\cA')
}$$
where all arrows except possibly $\sC_{\cA,\cA}$ are
equivalences; then the same must hold for $\sC_{\cA,\cA'}$.
\end{pfclaim}

\begin{claim}\label{cl_plus-commutes}
For every fibration $\cA\to\cC$ the diagram of functors :
$$
\xymatrix{ \sCart_\cC(\cE^+,\cA) \ar[rrr]^-{\sCart_\cC(j_\cE,\cA)}
\ar[d]_{\sC_{\cE^+,\cA}} & & & \sCart_\cC(\sC(\cE),\cA) \ar[d]^{(-)^+} \\
\sCart_\cC(\sC(\cE^+),\sC(\cA))
\ar[rrr]^-{\sCart_\cC(\omega_\cE,j_\cA)} & & &
\sCart_\cC(\sC(\cE)^+,\cA^+)
}$$
is essentially commutative (notation of claim
\ref{cl_again-plus-commutes}(i)).
\end{claim}
\begin{pfclaim} Let $G:\cE^+\to\cA$ be any $\cC$-cartesian
functor; the $\cC$-cartesian functor
$G':=(G\circ j_\cE)^+:\sC(\cE)^+\to\cA^+$ can be described as
follows. Let $[X,\psi:\cS\to\sC(\cE)]$ be any object of
$\sC(\cE)^+$; then
$G'([X,\psi])=[X,G\circ j_\cE\circ\psi:\cS\to\cA]$, and
notice that $j_\cE\circ\psi$ is the $\cC$-cartesian
functor given by the rule :
$$
(f:Y\to X)\mapsto[Y,\psi_f:\cC/Y\to\cE]
\qquad
\text{for every $f\in\Ob(\cS)$}.
$$
On the other hand, the $\cC$-cartesian functor
$G'':j_\cA\circ\sC(G)\circ\omega_\cE:\sC(\cE)^+\to\cA^+$
is described as follows. For any $[X,\psi]$ as in the
foregoing, $G''([X,\psi])=[X,G\circ\psi':\cC/X\to\cA]$,
where $\psi':\cS\to\cE^+$ is the $\cC$-cartesian functor
given by the rule :
$$
(f:Y\to X)\mapsto[Y,\sev^\cE\circ\psi\circ f_*]
\qquad
\text{for every $f\in\Ob(\cC/X)$}
$$
and notice that $\sev^\cE\circ\psi\circ f_*:\cS\times_Xf\to\cE$
is in turn given by the rule :
$$
(g:Z\to Y)\mapsto\psi_{f\circ g}(\one_Z)
\qquad
\text{for every $g\in\Ob(\cS\times_Xf)$}.
$$
Let $\psi'':\cS\to\cA$ be the restriction of $\psi'$; then
$[X,\psi']=[X,\psi'']$ in $\cE^+$, and $\cS\times_Xf=\cC/Y$
for every $(f:Y\to X)\in\Ob(\cS)$. For such $f$, we need to
compare the cartesian functors
$\psi_f,\sev^\cE\circ\psi\circ f_*:\cC/Y\to\cE$. However,
$\psi_f(\one_Y)=\sev^\cE\circ\psi\circ f_*(\one_Y)$, so there
exists a unique isomorphism of functors
$$
\tau_{\psi,f}:\sev^\cE\circ\psi\circ f_*\isom\psi_f
\qquad\text{such that}\qquad
(\tau_{\psi,f})_{\one_Y}=\one_{\psi_f(\one_Y)}.
$$
Explicitly, every object $(g:Z\to Y)$ of $\cC/Y$ yields a
morphism $g/Y:g\to\one_Y$ in $\cC/Y$, and therefore
$(\tau_{\psi,f})_g$ is determined by the commutativity of the
diagram :
$$
\xymatrix{
\psi_{f\circ g}(\one_Z) \ar[rr]^-{(\tau_{\psi,f})_g}
\ar[d]_{(\sev^\cE\circ\psi*f_*)_{g/Y}} & & \psi_f(g) \ar[d]^{\psi_f(g/Y)} \\
\psi_f(\one_Y) \rrdouble & & \psi_f(\one_Y)
}$$
However,
$(\sev^\cE\circ\psi\circ f_*)_{g/Y}=\psi_f(g/Y)\circ\psi_{g/X,\one_Z}$,
so that
$$
(\tau_{\psi,f})_g=\psi_{g/X,\one_Z}
\qquad
\text{for every $(g:Z\to Y)\in\Ob(\cC/Y)$}.
$$
We claim that we get an isomorphism of $\cC$-cartesian functors
$$
G([\tau_{\psi}]):G''([X,\psi])\isom G'([X,\psi])
\quad
f\mapsto G([\one_Y,\tau_{\psi,f}])
\quad
\text{for every $(f:Y\to X)\in\Ob(\cS)$}.
$$
Indeed, let $g/X:(f':Y'\to X)\to(f:Y\to X)$ be any morphism in
$\cS$; the assertion comes down to the commutativity of the
diagram :
$$
\xymatrix{ G([Y',\sev^\cE\circ\psi\circ f'_*]) \rrdouble
\ar[d]_{G([\one_{Y'},\tau_{\psi,f'}])} & &
G([Y',\sev^\cE\circ\psi\circ f_*\circ g_*])
\ar[d]^{G([\one_{Y'},\tau_{\psi,f}*g_*])} \\
G([Y',\psi_{f'}]) \ar[rr]^-{G([\one_{Y'},\psi_{g/X}])} & &
G([Y',\psi_f\circ g_*]).
}$$
The latter in turn follows by applying $G$ to the system
of identities :
$$
\psi_{g/X,h}\circ\psi_{h/X,\one_Z}=\psi_{g\circ h/X,\one_Z}
\qquad
\text{for every $(h:Z\to Y')\in\Ob(\cC/Y')$}
$$
which hold by functoriality of $\psi$, applied to the
morphisms $f'\circ h\xrightarrow{h/X}f'\xrightarrow{g/X}f$
of $\cS$.

Next, we check the naturality of the rule
$\psi\mapsto G([\tau_\psi])$. Hence, let
$[X_i,\psi_i:\cS_i\to\sC(\cE)]$ for $i=1,2$ be two objects of
$\sC(\cE)^+$, and $[t,\sigma]:[X_2,\psi_2]\to[X_1,\psi_1]$ a
morphism; {\em i.e.} $t:X_2\to X_1$ is a morphism of $\cC$
and $\sigma:\psi_{2|\cS_3}\Rightarrow(\psi_1\circ t_*)_{|\cS_3}$
is a natural $\cC$-transformation defined on a sieve
$\cS_3\subset\cS_2\cap(\cS_1\times_{X_1}t)$. We have
$$
G'([t,\sigma])=[t,(G\circ j_\cE)*\sigma:
G\circ j_\cE\circ\psi_{2|\cS_3}\Rightarrow
G\circ j_\cE\circ(\psi_1\circ t_*)_{|\cS_3}].
$$
On the other hand, notice that $\omega_\cE([t,\sigma]):
\omega_\cE([X_2,\psi_2])\to\omega_\cE([X_1,\psi_1])$ is the
natural transformation $\sigma':\psi'_2\Rightarrow\psi'_1*t_*$,
where $\psi'_i$ is the functor such that
$(f:Y\to X_i)\mapsto[Y,\sev^\cE\circ\psi_i\circ f_*]$ for
every $f\in\Ob(\cC/X_i)$, and $\sigma'_f:=\sev^\cE*\sigma*f_*$,
with $f_*:\cS_3\times_{X_2}f\to\cS_3$ for every $f\in\Ob(\cC/X_2)$.
Explicitly, $\sigma'_f$ assigns to every
$(g:Z\to Y)\in\Ob(\cS_3\times_{X_2}f)$ the morphism
$$
\sigma_{f\circ g,\one_Z}:
\psi_{2,f\circ g}(\one_Z)\to\psi_{1,t\circ f\circ g}(\one_Z).
$$
Thus, we need to check the commutativity of the diagram :
$$
\xymatrix{ [X_2,G\circ\psi'_2] \ar[rr]^-{G([\tau_{\psi_2}])}
\ar[d]_{[t,G*\sigma']} & &
[X_2,G\circ j_\cE\circ\psi_2] \ar[d]^{[t,(G\circ j_\cE)*\sigma]} \\
[X_1,G\circ\psi'_1] \ar[rr]^-{G([\tau_{\psi_1}])} & &
[X_1,G\circ j_\cE\circ\psi_1].
}$$
which in turn follows from the commutativity of the diagram :
$$
\xymatrix{ \psi_{2,f\circ g}(\one_Z) \ar[rr]^-{\sigma_{f\circ g,\one_Z}}
\ar[d]_{\psi_{2,g/X_2,\one_Z}} & & \psi_{1,t\circ f\circ g}(\one_Z)
\ar[d]^{\psi_{1,g/X_1,\one_Z}} \\
\psi_{2,f}(g) \ar[rr]^-{\sigma_{f,g}} & & \psi_{1,t\circ f}(g)
}$$
for every $f\in\Ob(\cC/X_2)$ and every $g\in\Ob(\cS_3\times_{X_2}f)$.
Since $g/X_1=t_*(g/X_2)$ in $\cS_1$, the assertion holds by
naturality of $\sigma$. Lastly, the naturality with respect
to natural $\cC$-transformations $G\Rightarrow G'$ is immediate
from the definitions.
\end{pfclaim}

Now, let $\cA\to\cC$ be a stack; then $j_\cA$ is an equivalence,
by lemma \ref{lem_bootstrap}(iv), and it follows easily that
$\sCart_\cC(\omega_\cE,j_\cA)$ is an equivalence as well.
The same holds for $\sC_{\cE^+,\cA}$ and $\sC_{\sC(\cE),\cA}$,
by claim \ref{cl_again-plus-commutes}(i). We deduce that
the functor $(-)^+$ appearing in the diagrams of claims
\ref{cl_plus-commutes} and \ref{cl_again-plus-commutes}(ii)
admits both left and right quasi-inverse functors, so it
is an equivalence. By claim \ref{cl_plus-commutes}, we finally
conclude that $\sCart_\cC(j_\cE,\cA)$ is an equivalence as well.
Now, set
$$
\cE^a:=\cE^{+++}.
$$
The composition $\sCart_\cC(j_\cE,\cA)\circ
\sCart_\cC(j_{\cE^+},\cA)\circ\sCart_\cC(j_{\cE^{++}},\cA)$ is
an equivalence
$$
\lambda_{\cE,\cA}:
\sCart_\cC(\cE^a,\cA)\isom\sCart_\cC(\sC(\sC(\sC(\cE))),\cA).
$$
The rules $(\cE,\cA)\mapsto\sCart_\cC(\cE^a,\cA)$ and
$(\cE,\cA)\mapsto\sCart_\cC(\sC(\sC(\sC(\cE))),\cA)$ yield
strict pseudo-functors
$$
\sCart_\cC((-)^a,-),\sCart_\cC(\sC\circ\sC\circ\sC,-):
\Fib(\cC)^o\times\Stack(\cC,J)\to\bCat
$$
and clearly the rule $(\cE,\cA)\mapsto\lambda_{\cE,\cA}$ defines
a strict pseudo-natural equivalence of functors
$$
\lambda:\sCart_\cC((-)^a,-)\Rightarrow\sCart_\cC(\sC\circ\sC\circ\sC,-).
$$
On the other hand, by composing evaluation functors we get
as well a strict pseudo-natural equivalence
$\sC\circ\sC\circ\sC\Rightarrow\one_{\Fib(\cC)}$, whence an
induced strict pseudo-natural equivalence
$$
\lambda':\sCart_\cC(-,\sF)\isom\sCart_\cC(\sC\circ\sC\circ\sC,-)
$$
where $\sF:\Stack(\cC,J)\to\Fib(\cC)$ is the (forgetful)
inclusion strict pseudo-functor. After choosing a quasi-inverse
for $\lambda$ (see theorem \ref{th_pseudo-nat-equiv}) and
composing with $\lambda'$, we obtain the required (non-strict)
pseudo-natural equivalence
$$
\sCart_\cC(-,\sF)\isom\sCart_\cC((-)^a,-).
$$
Lastly, the essential commutativity of \eqref{eq_free-palestine}
follows directly from lemma \ref{lem_bootstrap}(v).
\end{proof}

\subsection{Covering morphisms of prestacks}
In this section we wish to characterize the $i$-faithful functors
of prestacks (for $i=0,1,2$), by means of certain conditions that
are analogous to those given in definition \ref{def_tra-la-la} for
presheaves; to this aim, we make the following :

\begin{definition}\label{def_i-coverings}
Let $(\cC,J)$ be a site, $F:\cE\to\cC$ and $F':\cE'\to\cC$
two fibrations, and $\phi:\cE'\to\cE$ a $\cC$-cartesian functor.
\begin{enumerate}
\item
$\phi$ is a {\em $0$-covering} if for every $Y\in\Ob(\cE)$ there
exist a covering family $(f_i:X_i\to FY~|~i\in I)$ and for every
$i\in I$ an object $Y'_i\in\Ob(\cE'_{X_i})$ and a $\cC$-cartesian
morphism $h_i:\phi Y'_i\to Y$ in $\cE$ such that $F(h_i)=f_i$.
\item
$\phi$ is a {\em $1$-covering} if for every $X\in\Ob(\cC)$,
every $Y_1,Y_2\in\Ob(\cE'_X)$ and every morphism
$g:\phi Y_1\to\phi Y_2$ in $\cE_X$ there exist a covering
family $(f_i:X_i\to X~|~i\in I)$, and morphisms
$Y_2\xleftarrow{g_i}Z_i\xrightarrow{h_i}Y_1$ in $\cE'$,
where $h_i$ is $\cC$-cartesian for every $i\in I$, such that
$$
\phi(g_i)=g\circ\phi(h_i)
\qquad\text{and}\qquad
F'(h_i)=f_i.
$$
\item
$\phi$ is a {\em $2$-covering} if for every $X\in\Ob(\cC)$,
every $Y_1,Y_2\in\Ob(\cE'_X)$, and every pair $g_1,g_2:Y_1\to Y_2$
of morphisms in $\cE'_X$ with $\phi(g_1)=\phi(g_2)$ there
exist a covering family $(f_i:X_i\to X~|~i\in I)$, and for
every $i\in I$ a $\cC$-cartesian morphism
$$
h_i:Z_i\to Y_1
\quad\text{in $\cE'$ such that}
\quad
g_1\circ h_i=g_2\circ h_i
\quad\text{and}\quad
F'(h_i)=f_i.
$$
\end{enumerate}
\end{definition}

\begin{remark}\label{rem_realign}
Suppose that $(\cE\xrightarrow{F}\cC,\blambda)$ and
$(\cE'\xrightarrow{F'}\cC,\blambda')$ are split
fibrations over $\cC$, and $\phi:(\cE',\blambda')\to(\cE,\blambda)$
is a split cartesian functor. Let $\sc$ and $\sc'$ be the strict
pseudo-functors associated with $\blambda$ and $\blambda'$. Then
we can rephrase the conditions of definition \ref{def_i-coverings}
as follows :
\begin{enumerate}
\item
$\phi$ is $0$-covering if and only if for every $Y\in\Ob(\cE)$
there exist a covering family $(f_i:X_i\to FY~|~i\in I)$ and for
every $i\in I$ an object $Y'_i\in\Ob(\cE'_{X_i})$ and an isomorphism
$\phi Y'_i\isom\sc_{f_i}Y$ in $\cE_{X_i}$.
\item
$\phi$ is $1$-covering if and only if for every $X\in\Ob(\cC)$,
every $Y_1,Y_2\in\Ob(\cE'_X)$ and every morphism
$g:\phi Y_1\to\phi Y_2$ in $\cE_X$ we have a covering
family $(f_i:X_i\to X~|~i\in I)$, and for every $i\in I$ a
morphism $g_i:\sc'_{f_i}Y_1\to\sc'_{f_i}Y_2$ in $\cE'_{X_i}$ such
that $\phi(g_i)=\sc_{f_i}(g)$.
\item
$\phi$ is a {\em $2$-covering} if for every $X\in\Ob(\cC)$,
every $Y_1,Y_2\in\Ob(\cE'_X)$, and every pair $g_1,g_2:Y_1\to Y_2$
of morphisms in $\cE'_X$ with $\phi(g_1)=\phi(g_2)$ there
exists a covering family $(f_i:X_i\to X~|~i\in I)$ such that
$\sc'_{f_i}(g_1)=\sc'_{f_i}(g_2)$.
\end{enumerate}
\end{remark}

\begin{lemma}\label{lem_reduce-to-split-fibs}
Consider a site $(\cC,J)$, and an essentially commutative
diagram of the $2$-cate\-gory $\Fib(\cC)$, whose vertical
arrows are equivalences of categories :
$$
\xymatrix{ \cE' \ar[r]^-\phi \ar[d]_{\omega'} &
\cE \ar[d]^\omega \\
\cF' \ar[r]^-\psi & \cF.
}$$ 
Let also $i\in\{0,1,2\}$. Then $\phi$ is $i$-covering if and
only if the same holds for $\psi$.
\end{lemma}
\begin{proof} Denote $\cE\xrightarrow{F}\cC\xleftarrow{F'}\cE'$
and $\cF\xrightarrow{G}\cC\xleftarrow{G'}\cF'$ the respective
fibrations. The assumption means that there exists an isomorphism
of functors $\beta:\psi\circ\omega'\isom\omega\circ\phi$ such
that $G*\beta=\one_{F'}$.

Suppose that $\phi$ is $0$-covering, let $Y\in\Ob(\cF)$ be
any object, and set $X:=GY$; by assumption there exists
$Z\in\Ob(\cE_X)$ with an isomorphism $t:\omega Z\isom Y$ in
$\cF_X$ (corollary \ref{cor_fibrations}(i)). Then there exist
a covering family $(f_j:X_j\to X~|~j\in I)$ and for every
$j\in I$ a cartesian morphism $h_j:\phi Z'_j\to Z$ with $Fh_j=f_j$.
It follows that
$$
h'_j:=t\circ\omega(h_j)\circ\beta_{Z'_j}:\psi(\omega'Z'_j)\to Y
$$
is cartesian and $Gh'_j=f_j$ for every $j\in I$, so $\psi$ is
$0$-covering.

Conversely, suppose that $\psi$ is $0$-covering; let
$Y\in\Ob(\cE)$ and $X:=FY$. By assumption there exist
a covering family $(f_j:X_j\to X~|~j\in I)$ and for every
$j\in I$ a cartesian morphism $h_j:\psi Z_j\to\omega Y$ with
$Gh_j=f_j$. Then we find for every $j\in I$ an isomorphism
$t_j:\omega'Z'_j\isom Z_j$ in $\cF'_{X_j}$; it follows that
$h_j:=h\circ\psi(t_j)\circ\beta^{-1}_{Z'_j}:
\omega(\phi Z'_j)\to\omega Y$ is cartesian, and since
$\omega$ is an equivalence, there exists a unique cartesian
morphism $h'_j:\phi Z'_j\to Y$ in $\cE$ with $\omega(h'_j)=h_j$,
for every $j\in I$ (details left to the reader); by construction,
$Fh'_j=f_j$, so $\phi$ is $0$-covering.

Next, if $\phi$ is $1$-covering, let $X\in\Ob(\cC)$ and
$g:\psi Y_1\to\psi Y_2$ any morphism in $\cF_X$. We find
isomorphisms $t_i:\omega'Z_i\isom Y_i$ in $\cF'_X$ ($i=1,2$),
and a morphism $g':\psi(\omega'Z_1)\to\psi(\omega'Z_2)$
such that $\psi(t_2)\circ g'=g\circ\psi(t_1)$; then
$\beta_{Z_2}\circ g'=\omega(g'')\circ\beta_{Z_1}$ for a
morphism $g'':\phi Z_1\to\phi Z_2$ in $\cE_X$. By assumption,
there exist a covering family $(f_j:X_j\to X~|~j\in I)$ and
morphisms $Z_2\xleftarrow{g_j}W_j\xrightarrow{h_j}Z_1$ in $\cE'$
for every $j\in I$, such that $h_j$ is cartesian with $F'h_j=f_j$,
and $\phi(g_j)=g''\circ\phi(h_j)$. Set $h'_j:=t_1\circ\omega'(h_j)$
and $g'_j:=t_2\circ\omega'(g_j)$; then $h'_j$ is cartesian with
$G'(h'_j)=f_j$, and $\psi(g'_j)=g\circ\psi(h'_j)$ for every $j\in I$.
This shows that $\psi$ is $1$-covering.

Conversely, suppose that $\psi$ is $1$-covering, and let
$g:\phi Y_1\to\phi Y_2$ be a morphism in $\cE_X$. Then there
exist a covering family $(f_j:X_j\to X~|~j\in I)$ and morphisms
$\omega' Y_2\xleftarrow{g_j}Z_j\xrightarrow{h_j}\omega'Y_1$ in
$\cF'$ such that $h_j$ is cartesian with $G'h_j=f_j$ and
$\beta_{Y_2}\circ\psi(g_j)=\omega(g)\circ\beta_{Y_1}\circ\psi(h_j)$
for every $j\in I$. We pick also an isomorphism
$t_j:\omega'Z'_j\isom Z_j$ in $\cF'_{X_j}$ for
every $j\in I$; then there exist a unique cartesian morphism
$h'_j:Z'_j\to Y_1$ with $\omega'(h'_j)=h_j\circ t_j$ and a
unique morphism $g'_j:Z'_j\to Y_2$ with $\omega'(g'_j)=g_j\circ t_j$.
A direct computation shows that
$\omega(g\circ\phi(h'_j))=\omega(\phi(g'_j))$, whence
$g\circ\phi(h'_j)=\phi(g'_j)$ for every $j\in I$. Thus,
$\phi$ is $1$-cartesian.

Lastly, suppose that $\phi$ is $2$-covering, and let
$g_1,g_2:Y_1\to Y_2$ be two morphisms in $\cF'_X$ such
that $\psi(g_1)=\psi(g_2)$. Pick isomorphisms
$t_i:\omega'Z_i\isom Y_i$ for $i=1,2$; then there exists
for $i=1,2$ a unique morphism $g'_i:Z_1\to Z_2$ such that
$t_2\circ\omega'(g'_i)=g_i\circ t_1$. It follows that
$$
\begin{aligned}
\psi(t_2)\circ\beta^{-1}_{Z_2}\circ\omega(\phi(g'_1))
&\,=\psi(t_2)\circ\psi(\omega'g_1)\circ\beta^{-1}_{Z_1} \\
&\,=\psi(g_1\circ t_1)\circ\beta^{-1}_{Z_1} \\
&\,=\psi(g_2\circ t_1)\circ\beta^{-1}_{Z_1} \\
&\,=\psi(t_2)\circ\beta^{-1}_{Z_2}\circ\omega(\phi(g'_2))
\end{aligned}
$$
whence $\phi(g'_1)=\phi(g'_2)$. Then we have a covering family
$(f_j:X_j\to X~|~j\in I)$ and for every $j\in I$ a cartesian
morphism $h_j:W_j\to Z_1$ such that $g'_1\circ h_j=g'_2\circ h_j$
and $Fh_j=f_j$. Hence $h'_j:=t_1\circ\omega'(h_j):\omega'W_j\to Y_1$
is cartesian with $G'h'_j=f_j$ and $g_1\circ h'_j=g_2\circ h'_j$.
This proves that $\psi$ is $2$-covering. Conversely, suppose
that $\psi$ is $2$-covering, and let $g_1,g_2:Y_1\to Y_2$ be
two morphisms in $\cE'_X$ such that $\phi(g_1)=\phi(g_2)$.
Set $g'_i:=\omega'(g_i)$ for $i=1,2$; then $\psi(g'_1)=\psi(g'_2)$,
so there exist a covering family $(f_j:X_j\to X~|~j\in I)$
and for every $j\in I$ a cartesian morphism $h_j:Z_j\to\omega' Y_1$
such that $g'_1\circ h_j=g'_2\circ h_j$ and $G'h_j=f_j$. For
every $j\in I$ pick an isomorphism $t_j:\omega' Z'_j\isom Z_j$
in $G'^{-1}X$; there exists a unique cartesian morphism
$h'_j:Z'_j\to Y_1$ such that $\omega'(h'_j)=h_j\circ t_j$.
It is easily seen that
$\omega'(g_1\circ h'_j)=\omega'(g_2\circ h'_j)$, whence
$g_1\circ h'_j=g_2\circ h'_j$, and $F'h'_j=f_j$ for every
$j\in I$. Thus, $\phi$ is $2$-covering.
\end{proof}

\begin{lemma}\label{lem_yoga-i-coverings}
Let $(\cC,J)$ be a site,
$F:\cE\to\cC,F':\cE'\to\cC,F'':\cE''\to\cC$ three
fibrations, $\phi:\cE''\to\cE'$ and $\psi:\cE'\to\cE$
two cartesian functors, and $i\in\{0,1,2\}$. We have:
\begin{enumerate}
\item
If $\phi$ and $\psi$ are $i$-covering, the same holds for
$\psi\circ\phi$.
\item
If $\psi\circ\phi$ is $i$-covering, and $\phi$ is $j$-covering
for every $j<i$, then $\psi$ is $i$-covering.
\item
If $\psi\circ\phi$ is $i$-covering and $\psi$ is $j$-covering
for every $j>i$, then $\phi$ is $i$-covering.
\end{enumerate}
\end{lemma}
\begin{proof} In view of lemma \ref{lem_reduce-to-split-fibs}
and claim \ref{cl_split-fibrations} we may replace
$\phi$ and $\psi$ by $\sC(\phi)$ and $\sC(\psi)$, and
assume from start that $\phi$ and $\psi$ are split
cartesian functors. We then denote by $\sc$, $\sc'$ and
$\sc''$ the strict pseudo-functors associated with the
split cleavages of $\cE$, $\cE'$ and $\cE''$ respectively.

(i): Suppose first that $i=0$, let $Y\in\Ob(\cE)$ be any
object, and set $X:=FY$; by assumption there exist a covering
family $(f_j:X_j\to X~|~j\in I)$ and for every $j\in I$ an
isomorphism $h_j:\psi Y'_j\isom\sc_{f_j}Y$ in $\cE_{X_j}$.
Likewise, for every $j\in I$ there exist a covering family
$(f_{j\lambda}:X_{j\lambda}\to X_j~|~\lambda\in\Lambda_j)$ and for
every $\lambda\in\Lambda_j$ an isomorphism
$h_{j\lambda}:\phi Y''_{j\lambda}\isom\sc'_{f_{j\lambda}}Y'_j$ in
$\cE'_{X_{j\lambda}}$. Then the family
$(f'_{j\lambda}:=f_j\circ f_{j\lambda}~|~j\in I,\ \lambda\in\Lambda_j)$
covers $X$ and $\sc_{f_{j\lambda}}(h_j)\circ\psi(h_{j\lambda}):
\psi\circ\phi(Y''_{j\lambda})\to\sc_{f'_{j\lambda}}Y$ is an isomorphism
in $\cE_{X_{j\lambda}}$ for every $j\in I$ and $\lambda\in\Lambda_j$,
whence the contention, in this case.

Next, suppose that $i=1$. Let $X\in\Ob(\cC)$,
$Y''_1,Y''_2\in\Ob(\cE''_X)$, and
$g:\psi\circ\phi Y''_1\to\psi\circ\phi Y''_2$ any morphism
in $\cE_X$. By assumption there exist a covering family
$(f_j:X_j\to X~|~j\in I)$ and for every $j\in I$ a morphism
$$
g_j:\phi(\sc''_{f_j}Y''_1)=\sc'_{f_j}(\phi Y''_1)\to
\sc'_{f_j}(\phi Y''_2)=\phi(\sc''_{f_j}Y''_2)
\quad\text{in $\cE'_{X_j}$ such that}\quad
\psi(g_j)=\sc_{f_j}(g).
$$
Since $\phi$ is $1$-covering, we find for every $j\in I$ a
covering family
$(f_{j\lambda}:X_{j\lambda}\to X_j~|~\lambda\in\Lambda_j)$ and
for every $\lambda\in\Lambda_j$ a morphism
$$
g_{j\lambda}:\sc''_{f_{j\lambda}}(\sc''_{f_j}Y''_1)\to
\sc''_{f_{j\lambda}}(\sc''_{f_j}Y''_2)
\qquad\text{in $\cE'_{X_j}$ such that}\qquad
\phi(g_{j\lambda})=\sc'_{f_{j\lambda}}(g_j).
$$
Then the family
$(f'_{j\lambda}:=f_j\circ f_{j\lambda}~|~j\in I,\ \lambda\in\Lambda_j)$
covers $X$, and $\psi\circ\phi(g_{j\lambda})=\sc_{f'_{j\lambda}}(g)$
for every $j\in I$ and every $\lambda\in\Lambda_j$. It follows
that $\psi\circ\phi$ is $1$-covering.

Lastly, let $i=2$; consider morphisms $g_1,g_2:Y''_1\to Y''_2$
in $\cE''_X$ with $\psi\circ\phi(g_1)=\psi\circ\phi(g_2)$. By
assumption, there exists a covering family $(f_j:X_j\to X~|~j\in I)$
such that
$$
\phi(\sc''_{f_j}(g_1))=\sc'_{f_j}(\phi(g_1))=\sc'_{f_j}(\phi(g_2))
=\phi(\sc''_{f_j}(g_2))
\qquad
\text{for every $j\in I$}.
$$
We may then find for every $j\in I$ a covering family
$(f_{j\lambda}:X_{j\lambda}\to X_j~|~\lambda\in\Lambda_j)$ such that
$\sc''_{f_{j\lambda}}\circ\sc''_{f_j}(g_1)=
\sc''_{f_{j\lambda}}\circ\sc''_{f_j}(g_2)$ for every $\lambda\in\Lambda_j$.
We conclude that $\psi\circ\phi$ is $2$-covering.

(ii): Again, we suppose first that $i=0$, in which case
the assumption means that $\psi\circ\phi$ is $0$-covering
(and there is no condition on $\phi$), and we need to check
that $\psi$ is $0$-covering. However, the assertion follows
straightforwardly from the definitions (details left to the reader).

Next, suppose that $i=1$, in which case $\psi\circ\phi$ is
$1$-covering and $\phi$ is $0$-covering. Let
$g:\psi Y'_1\to\psi Y'_2$ be a morphism in $\cE_X$, for some
$X\in\Ob(\cC)$; by assumption there exist for $t=1,2$ a
covering family $(f_{tj}:X_{tj}\to X~|~j\in I_t)$ and
isomorphisms $h_{tj}:\phi Y''_{tj}\isom\sc'_{f_j}Y'_t$ in
$F'^{-1}X_j$ for every $j\in I_t$. Notice then that for every
morphism $f':X'\to X_{ti}$ in $\cC$ we get an isomorphism
$\sc'_{f'}(h_{tj}):\phi(\sc''_{f'}Y''_{tj})\isom\sc'_{f_j\circ f'}Y'_t$
in $\cE'_{X'}$. We deduce that, after replacing the covering
families $f_{1\bullet}$ and $f_{2\bullet}$ by a common refinement,
we may assume that $I:=I_1=I_2$ and $f_j:=f_{1j}=f_{2j}$ for
every $j\in I$. In this situation, for every $j\in I$ there
exists a morphism
$$
g_j:\psi\circ\phi(Y''_{1j})\to\psi\circ\phi(Y''_{2j})
\qquad\text{in $\cE_{X_j}$ such that}\qquad
\psi(h_{2j})\circ g_j=\sc_{f_j}(g)\circ\psi(h_{1j}).
$$
By our assumption on $\psi\circ\phi$, for every $j\in I$
we then find a covering family
$(f_{j\lambda}:X_{j\lambda}\to X_j~|~\lambda\in\Lambda_j)$ and for
every $\lambda\in\Lambda_j$ a morphism
$$
g_{j\lambda}:Y''_{1j}\to Y''_{1j}
\qquad\text{in $\cE''_{X_{j\lambda}}$}\qquad\text{such that}\qquad
\psi\circ\phi(g_{j\lambda})=\sc_{f_{j\lambda}}(g_j).
$$
It follows that $\psi(\sc'_{f_{j\lambda}}h_{2j})\circ\psi(\phi(g_{j\lambda}))=
\sc_{f_j\circ f_{j\lambda}}(g)\circ\psi(\sc'_{f_j}h_{1j})$ for every
$j\in I$ and $\lambda\in\Lambda_j$. Since the family
$(f_i\circ f_{i\lambda}~|~i\in I,\ \lambda\in\Lambda_i)$
covers $X$, we conclude that $\psi$ is $1$-covering.

Next, suppose that $i=2$, so $\psi\circ\phi$ is $2$-covering
and $\phi$ is $j$-covering for $j=0,1$. Let $X\in\Ob(\cC)$ and
$g_1,g_2:Y'_1\to Y'_2$ two morphisms of $\cE'_X$ such that
$\psi(g_1)=\psi(g_2)$. Arguing as in the previous case, we
find a covering family $(f_j:X_j\to X~|~j\in I)$ and isomorphisms
$h_{tj}:\phi Y''_{tj}\isom\sc'_{f_j}Y'_t$ in $\cE'_{X_j}$  for $t=1,2$
and every $j\in I$. In this situation, for every $j\in I$ and
$t=1,2$ there exists a unique morphism
$g_{tj}:\phi Y''_{1i}\to\phi Y''_{2i}$ in $\cE'_{X_j}$ such that
\set\begin{equation}\label{eq_dont-want-to-fuck-you}
h_{2j}\circ g_{tj}=\sc_{f_j}(g_t)\circ h_{1j}.
\end{equation}
Notice that $\psi(h_{2j})\circ\psi(g_{1j})=\sc'_{f_j}(\psi(g_1))=
\sc'_{f_j}(\psi(g_2))=\psi(h_{2j})\circ\psi(g_{2j})$, whence
$$
\psi(g_{1j})=\psi(g_{2j})
\qquad
\text{for every $j\in I$}.
$$
Then, since $\phi$ is $1$-covering, for every $j\in I$ and
$t=1,2$ we find a covering family $f_{tj\bullet}:=
(f_{tj\lambda}:X_{tj\lambda}\to X_j~|~\lambda\in\Lambda_{tj})$ and
morphisms
$$
g_{tj\lambda}:\sc''_{f_{tj\lambda}}Y''_{1j}\to\sc''_{f_{tj\lambda}}Y''_{2j}
\qquad\text{in $\cE''_{X_{tj\lambda}}$}\qquad\text{such that}\qquad
\phi(g_{tj\lambda})=\sc'_{f_{tj\lambda}}(g_{tj})
$$
and arguing as in the foregoing, we may replace $f_{1j\bullet}$
and $f_{2j\bullet}$ by a common refinement, and assume that
$\Lambda_j:=\Lambda_{1j}=\Lambda_{2j}$ for every $j\in I$, and
$f_{j\lambda}:=f_{1j\lambda}=f_{2j\lambda}$, for every $j\in I$ and
every $\lambda\in\Lambda_j$. We deduce that
$$
\psi\circ\phi(g_{1j\lambda})=\sc''_{f_{j\lambda}}(\psi(g_{1j}))=
\sc''_{f_{j\lambda}}(\psi(g_{2j}))=\psi\circ\phi(g_{2j\lambda})
\qquad
\text{for every $j\in I$ and $\Lambda\in\Lambda_j$}.
$$
Since $\psi\circ\phi$ is $2$-covering, there exists therefore
a covering family
$(f_{j\lambda\lambda'}:X_{j\lambda\lambda'}\to X_{j\lambda}~|~
\lambda'\in\Lambda_{j\lambda})$ such that
$\sc''_{f_{j\lambda\lambda'}}(g_{1j\lambda})=
\sc''_{f_{j\lambda\lambda'}}(g_{2j\lambda})$ for every
$\lambda'\in\Lambda_{j\lambda}$. Therefore
$\sc'_{f_{j\lambda\lambda'}}\circ\sc'_{f_{j\lambda}}(g_{1j})=
\phi(\sc''_{f_{j\lambda\lambda'}}(g_{1j\lambda}))=
\phi(\sc''_{f_{j\lambda\lambda'}}(g_{2j\lambda}))=
\sc'_{f_{j\lambda\lambda'}}\circ\sc'_{f_{j\lambda}}(g_{2j})$ for every
$j\in I$, $\lambda\in\Lambda_j$ and $\lambda'\in\Lambda_{j\lambda}$.
Set $f'_{j\lambda\lambda'}:=f_j\circ f_{j\lambda}\circ f_{j\lambda\lambda'}$
for every such $j,\lambda,\lambda'$; combining with
\eqref{eq_dont-want-to-fuck-you}, we deduce easily that
$\sc_{f'_{j\lambda\lambda'}}(g_1)=\sc_{f'_{j\lambda\lambda'}}(g_2)$, and
we conclude that $\psi$ is $2$-covering, as sought.

(iii): The case where $i=2$ is immediate from the definitions.

Say that $i=1$, so $\psi\circ\phi$ is $1$-covering and
$\psi$ is $2$-covering; let $X\in\Ob(\cC)$ and
$g:\phi Y''_1\to\phi Y''_2$ a morphism in $\cE'_X$. By
assumption, there exist a covering family $(f_j:X_j\to ~|~j\in I)$
and a morphism $g_j:\sc_{f_j}Y''_1\to\sc_{f_j}Y''_2$ in $\cE''_{X_j}$
for every $j\in I$ such that
$$
\psi\circ\phi(g_j)=\sc_{f_j}(\psi(g))=\psi(\sc'_{f_j}(g)).
$$
Then, since $\psi$ is $2$-covering, we find for every $j\in I$ a
covering family $(f_{j\lambda}:X_{j\lambda}\to X_j~|~\lambda\in\Lambda_j)$
such that $\sc'_{f_{j\lambda}}(\phi(g_j))=\sc'_{f_{j\lambda}}\circ\sc'_{f_j}(g)$
for every $\lambda\in\Lambda_j$, {\em i.e.}
$\phi(\sc''_{f_{j\lambda}}(g_j))=\sc'_{f_j\circ f_{j\lambda}}(g)$.
This shows that $\phi$ is $1$-covering.

In order to deal with the case where $i=0$, let us remark :

\begin{claim}\label{cl_locally-isomorphic}
Suppose $\psi$ is $j$-covering for $j=1,2$. Let $X\in\Ob(\cC)$
and $Z,Z'\in\Ob(\cE'_X)$ such that $\psi Z$ and $\psi Z'$ are
isomorphic in $\cE_X$. Then there exists a covering family
$(f_j:X_j\to X~|~j\in I)$ such that $\sc'_{f_j}Z$ and
$\sc'_{f_j}Z'$ are isomorphic in $\cE'_{X_j}$ for every $j\in I$. 
\end{claim}
\begin{pfclaim} Let $l:\psi Z\isom\psi Z'$ be an isomorphism
in $\cE_X$; since $\psi$ is $1$-covering, there exist a
covering family $(f'_j:X'_j\to X~|~j\in I)$ and for every
$j\in I$ a morphism $l_j:\sc'_{f'_j}Z\to\sc'_{f'_j}Z'$ 
such that $\psi(l_j)=\sc_{f'_j}(l)$. Likewise, there exist
a covering family $(f''_j:X'_j\to X~|~j\in I')$ and for every
$j\in I'$ a morphism $l'_j:\sc'_{f''_j}Z'\to\sc'_{f''_j}Z$ such
that $\psi(l'_i)=\sc_{f''_j}(l^{-1})$. After replacing
$(f'_j~|~j\in I)$ and $(f''_j~|~j\in I')$ by a common
refinement, we may then assume that $I=I'$ and
$f_j:=f'_j=f''_j$ for every $j\in I$ (details left to the reader).
It follows that
$\phi(k_j\circ l_j)=\sc_{f_j}(\one_Z)=\phi(\one_{\sc'_{f_j}Z})$
and $\phi(l_j\circ k_j)=\sc_{f_j}(\one_{Z'})=\phi(\one_{\sc'_{f_j}Z'})$
for every $j\in I$. Then, since $\psi$ is $2$-covering,
for every $j\in I$ there exists a covering family
$(f_{j\lambda}:X_{j\lambda}\to X_j~|~\lambda\in\Lambda_j)$ such that
$\sc'_{f_{j\lambda}}(k_j\circ l_j)=\sc'_{f_{j\lambda}}(\one_{\sc'_{f_j}Z})$;
likewise, for every $j\in I$ there exists a covering family
$(f'_{j\lambda}:X'_{j\lambda}\to X_j~|~\lambda\in\Lambda'_j)$ such that
$\sc'_{f'_{j\lambda}}(l_j\circ k_j)=\sc'_{f'_{j\lambda}}(\one_{\sc'_{f_j}Z'})$.
After replacing these families by a common refinement, we may
assume that $\Lambda_j=\Lambda'_j$ for every $j\in I$, and
$f_{j\lambda}=f'_{j\lambda}$ for every $j\in I$ and every
$\lambda\in\Lambda_j$. Then, set $h_{j\lambda}:=f_j\circ f_{j\lambda}$
for every $j\in I$ and every $\lambda\in\Lambda_j$; we deduce
easily that $\sc'_{f'_{j\lambda}}(l_j)$ is an isomorphism
$\sc'_{h_{j\lambda}}Z\isom\sc'_{h_{j\lambda}}Z'$, whence the claim.
\end{pfclaim}

Now, suppose that $\psi\circ\phi$ is $0$-covering and $\psi$
is $j$-covering for $j=1,2$, and let $X\in\Ob(\cC)$,
$Y'\in\Ob(\cE'_X)$; by assumption there exist a covering
family $(f_j:X_j\to X~|~j\in I)$ and for every $j\in I$ an
isomorphism
$\psi\circ\phi Y''_j\isom\sc_{f_j}\psi Y'=\psi(\sc'_{f_j}Y')$
in $\cE_X$. By claim \ref{cl_locally-isomorphic}, there
exists for every $j\in I$ a covering family
$(f_{j\lambda}:X_{j\lambda}\to X_j~|~j\in I)$ such that
$\sc'_{f_{j\lambda}}(\phi Y''_j)=\phi(\sc''_{f_{j\lambda}}Y''_j)$ is
isomorphic to $\sc'_{f_j\circ f_{j\lambda}}Y'$ in $\cE'_{X_i}$ for
every $\lambda\in\Lambda_j$. This shows that $\phi$ is
$0$-covering, and concludes the proof of the lemma.
\end{proof}

\begin{lemma}\label{lem_charact-1-2-coverings}
Let  $C:=(\cC,J)$ be a site,
$\cE\xrightarrow{F}\cC\xleftarrow{F'}\cE'$ two fibrations,
and $\phi:\cE\to\cE'$ a $\cC$-cartesian functor. The following
conditions are equivalent :
\begin{enumerate}
\item
$\phi$ is $1$-covering (resp. $2$-covering).
\item
For every $X\in\Ob(\cC)$ and every $\sigma,\sigma'\in\cE(X)$,
the induced morphism :
\set\begin{equation}\label{eq_associated-cart-sections}
\cCart(\sigma,\sigma')^a\to\cCart(\phi\circ\sigma,\phi\circ\sigma')^a
\end{equation}
of sheaves on $C/X$ is an epimorphism (resp. a monomorphism).
\end{enumerate}
\end{lemma}
\begin{proof} Recall that -- as described in \eqref{eq_F-pariah} --
for every $X,\sigma,\sigma'$ as in (ii) we have a morphism of
presheaves  $\phi^*_{\sigma,\sigma'}:\cCart(\sigma,\sigma')\to
\cCart(\phi\circ\sigma,\phi\circ\sigma')$ on $\cC/X$, and
\eqref{eq_associated-cart-sections} is the morphism of associated
sheaves $(\phi^*_{\sigma,\sigma'})^a$ on the site $C/X$.
Then the assertion follows easily by inspecting the definitions,
and taking into account corollary \ref{cor_bicover}(i,ii) and
remark \ref{rem_iprippi}(ii,iii).
\end{proof}

\begin{proposition}\label{prop_faith-separ-cover}
Let $(\cC,J)$ be a site, $\cE\xrightarrow{F}\cC\xleftarrow{F'}\cE'$
two fibrations, $\phi:\cE\to\cE'$ a $\cC$-cartesian functor,
and $i\in\{0,1,2\}$ such that $\cE$ is $i$-separated
and $\cE'$ is $(i-1)$-separated. Then $\phi$ is $i$-faithful
if and only if it is $j$-covering for every $j\geq 2-i$.
\end{proposition}
\begin{proof} We consider first the case where $i=0$,
so by assumption $\cE$ is a $0$-separated prestack (and
no conditions on $\cE'$); the assertion is that $\phi$
is faithful if and only if it is $2$-covering. However,
if $\phi$ is faithful, it is obviously $2$-covering.
Thus, suppose that $\phi$ is $2$-covering, and let
$X\in\Ob(\cC)$ and $g_1,g_2:Y\to Z$ two morphisms in
$\cE_X$ such that $\phi(g_1)=\phi(g_2)$; by assumption
there exist a covering family $(f_j:X_j\to X~|~j\in I)$
and for every $j\in I$ a cartesian morphism
$h_j:Y_j\to Y$ in $\cE$ such that $Fh_j=f_j$ and
$g_1\circ h_j=g_2\circ h_j$. We may then find two cartesian
sections $\psi_Y,\psi_Z\in\cE(X)$ such that
$\psi_Y(\one_X)=Y$ and $\psi_Z(\one_X)=Z$, and for $i=1,2$
there exists a unique natural $\cC$-transformation
$\sigma_i:\psi_Y\Rightarrow\psi_Z$ such that
$\sigma_{i,\one_X}=g_i$. It follows that
$(\sigma_1)_{|\cS}=(\sigma_2)_{|\cS}$, so $\sigma_1=\sigma_2$,
since $\cE$ is $0$-separated; hence $g_1=g_2$, so $\phi$ is
faithful.

In case $i=1$, the prestack $\cE$ is $1$-separated and $\cE'$
is $0$-separated, and we need to check that $\phi$ is fully
faithful if and only if it is $j$-covering for $j=1,2$.
Again, if $\phi$ is fully faithful, obviously it is $j$-covering
for $j=1,2$. Conversely, if the latter condition holds, by the
foregoing case we know already that $\phi$ is faithful, so it
remains only to check that $\phi$ is full. Thus, let
$X\in\Ob(\cC)$ and $g:\phi Y_1\to\phi Y_2$ any morphism in
$\cE'_X$; by assumption, there exist a covering family
$(f_j:X_j\to X~|~j\in I)$ and morphisms
$Y_2\xleftarrow{g_j}Z_j\xrightarrow{h_j}Y_1$ with $h_j$ cartesian,
such that $Fh_j=f_j$ and $\phi(g_j)=g\circ\phi(h_j)$ for every
$j\in I$. As in the foregoing, we pick for $i=1,2$ a cartesian
section $\psi_i\in\cE(X)$ such that $\psi_i(\one_X)=Y_i$,
and then there exists a unique natural $\cC$-transformation
$\sigma:\phi\circ\psi_1\Rightarrow\phi\circ\psi_2$ such that
$\sigma_{\one_X}=g$. Since $\psi_1(f_j/X):\psi_1(f_j)\to Y_1$
is cartesian, there exists moreover a unique isomorphism in
$F^{-1}X_j$
$$
h'_j:\psi_1(f_j)\isom Z_j
\qquad\text{such that}\qquad
h_j\circ h'_j=\psi_1(f_j/X).
$$
Likewise, since $\psi_2(f_j/X):\psi_2(f_j)\to Y_2$ is cartesian,
there exists a unique morphism in $F^{-1}X_j$
$$
g'_j:\psi_1(f_j)\to\psi_2(f_j)
\qquad\text{such that}\qquad
g_j\circ h'_j=\psi_2(f_j/X)\circ g'_j.
$$
We get therefore for every $j\in I$ a unique natural
$\cC$-transformation
$$
\sigma_j:\psi_1\circ f_{j*}\Rightarrow\psi_2\circ f_{j*}
\qquad\text{such that}\qquad
\sigma_{j,\one_{X_j}}=g'_j.
$$

\begin{claim}\label{cl_lift-sigma}
$\phi*\sigma_j=\sigma*f_j$ for every $j\in I$.
\end{claim}
\begin{pfclaim} It suffices to check that
$\phi(\sigma_{j,\one_{X_j}})=(\sigma*f_j)_{\one_{X_j}}$, {\em i.e.}
$\phi(g'_j)=\sigma_{f_j}$ for every $j\in I$. However :
$\phi(\psi_2(f_j/X))\circ\phi(g')=\phi(g_j)\circ\phi(h'_j)=
g\circ\phi(h_j)\circ\phi(h'_j)=g\circ\phi(\psi_1(f_j/X))=
\phi(\psi_2(f_j/X))\circ\sigma_{f_j}$, whence the assertion,
since $\phi(\psi_2(f_j/X))$ is cartesian.
\end{pfclaim}

Recall now that the proof of theorem \ref{th_split-fibration}
yields a natural equivalence $\eta_{\sC\cE)}:\sC(\cE)\isom\sC(\sC(\cE))$,
induced by the strict pseudo-natural equivalence $\beta^\blambda$
associated with the natural split cleavage $\blambda$ of
$\sC(\cE)$. We set $\psi_i^*:=\eta_{\sC(\cE)}(\psi_i)$ for $i=1,2$.
Explicitly, $\psi^*_i(f)=\psi_i\circ f_*$ for every
$f\in\Ob(\cC/X)$, and $\psi^*_{i,g/X}=\one_{\psi^*_i(f')}$ for
every morphism $g/X:f'\to f$ in $\cC/X$. Let $\cS\subset\cC/X$
be the sieve generated by $(f_j~|~j\in I)$; we construct as follows
a natural $\cC$-transformation
$$
\sigma^\dagger:\psi^*_{|\cS}\Rightarrow\psi^*_{|\cS}.
$$
For every $(f:X'\to X)\in\Ob(\cS)$, pick $j\in J$ such that
there exists a factorisation $f=f_j\circ f'$ for some morphism
$f':X\to X_j$ of $\cC$, and set $\sigma^\dagger_f:=\sigma_j*f'_*$.
Let us check that $\sigma^\dagger_f$ is independent of the choice
of factorisation : indeed, if $f=f''\circ f_k$ for some other
$k\in I$ and some morphism $f'':X'\to X_k$, we have
$\phi*(\sigma_j*f')=\sigma*f_j*f'=\sigma*f_k*f''=\phi*(\sigma_k*f'')$
by claim \ref{cl_lift-sigma}, whence the contention, since $\phi$
is faithful. Likewise, the naturality of the rule
$f\mapsto\sigma^\dagger_f$ can be checked after composition with
$\phi$, where it follows easily (details left to the reader).

Now, since $\cE$ is $1$-separated, the same holds for
$\sC(\cE)$, hence $\sigma^\dagger$ extends to a natural
$\cC$-transformation $\sigma^*:\psi^*_1\Rightarrow\psi^*_2$.
Then there exists a unique natural $\cC$-transformation
$\sigma':\psi_1\Rightarrow\psi_2$ such that
$\eta_{\sC(\cE)}(\sigma')=\sigma^*$. We claim that
$$
\phi*\sigma'=\sigma.
$$
Indeed, it suffices to show that
$\eta_{\sC(\cE)}(\phi*\sigma')=\eta_{\sC(\cE)}(\sigma)$, {\em i.e.}
that $\sC(\phi)*\sigma^*=\eta_{\sC(\cE)}(\sigma)$, and since $\cE'$
is $0$-separated, we are reduced to checking that
$(\sC(\phi)*\sigma^*)_{|\cS}=\eta_{\sC(\cE)}(\sigma)_{|\cS}$, {\em i.e.}
that $\sC(\phi)*\sigma^\dagger=\eta_{\sC(\cE)}(\sigma)_{|\cS}$. Then we
are further reduced to showing that
$\phi*\sigma^\dagger_{f_j}=\sigma*f_{j*}$ for every $j\in I$. But
the latter identities hold by claim \ref{cl_lift-sigma}. Lastly,
set $g':=\sigma'_{\one_X}$; we deduce that $\phi(g')=g$, and this
shows that $\phi$ is full, as required.

In case $i=2$, the prestack $\cE'$ is $1$-separated and $\cE$
is a stack, and we need to check that $\phi$ is an equivalence
if and only if it is $j$-covering for $j=0,1,2$. Again, if
$\phi$ is an equivalence, clearly it is $j$-covering for
all $j$; conversely, if this condition holds, by the foregoing
we know already that $\phi$ is fully faithful, so it remains
only to check that $\phi$ is essentially surjective. To this
aim, let $X\in\Ob(\cC)$ and $Y'\in\Ob(\cE'_X)$; by assumption,
there exist a covering family $(f_j:X_j\to X~|~j\in I)$ and
for every $j\in I$ a cartesian morphism $h_j:\phi Y_j\to Y'$
with $F'h_j=f_j$. Pick, for every $j\in I$ a cartesian section
$\psi_j\in\cE(X_j)$ with $\psi_j(\one_{X_j})=Y_j$. Now, for
every $j,k\in I$ consider the category
$\cC/X_{jk}:=\cC/X_j\times_{(f_{j*},f_{k*})}\cC/X_k$ (see
\eqref{subsec_step-one}). The objects of $\cC/X_{jk}$ are the
pairs of morphisms $X_j\xleftarrow{f'_j}X'\xrightarrow{f'_k}X_k$
such that $f':=f_j\circ f'_j=f_k\circ f'_k$. We have two natural
projections $\cC/X_j\xleftarrow{\pi^0_{jk}}\cC/X_{jk}
\xrightarrow{\pi^1_{jk}}\cC/X_k$ that send every such pair
$(f'_j,f'_k)$ to $f'_j$ and respectively $f'_k$.
For every such pair $(f'_j,f'_k)$ and $s=j,k$ we get cartesian
morphisms $\psi_s(f'_s/X_s):\psi_s(f'_s)\to Y_s$ such that
$F\psi_s(f'_s/X_s)=f'_s$. There follows a unique isomorphism
$$
\tau:\phi\circ\psi_j(f'_j)\isom\phi\circ\psi_k(f'_k)
\qquad\text{such that}\qquad
h_k\circ\phi(\psi_k(f'_k/X_k))\circ\tau=
h_j\circ\phi(\psi_j(f'_j/X_j))
$$
and since $\phi$ is fully faithful, we get a unique isomorphism
$$
\omega^{ij}_{(f'_j,f'_k)}:\psi_j(f'_j)\isom\psi_k(f'_k)
\qquad\text{such that}\qquad
\phi(\omega^{ij}_{(f'_j,f'_k)})=\tau
$$
and we claim that the rule $(f'_j,f'_k)\mapsto\omega^{ij}_{(f'_j,f'_k)}$
yields an isomorphism of functors
$$
\omega^{ij}:\psi_j\circ\pi^0_{jk}(f'_j,f'_k)\isom
\psi_k\circ\pi^1_{jk}(f'_j,f'_k).
$$
Indeed, let $(X_j\xleftarrow{f''_j}X''\xrightarrow{f''_k}X_k)$
be another object of $\cC/X_{jk}$; a morphism
$(f''_j,f''_k)\to(f'_j,f'_k)$ is a morphism $t:X''\to X'$ in
$\cC$ such that $f'_j\circ t=f''_j$ and $f'_k\circ t=f''_k$,
and the assertion amounts to the identity :
$$
A:=\psi_k(t/X_k)\circ\omega^{jk}_{(f''_j,f''_k)}=
B:=\omega^{jk}_{(f'_j,f'_k)}\circ\psi_j(t/X_j).
$$
However, we have :
$$
\begin{aligned}
h_k\circ(\phi\circ\psi_k(f'_k/X_k))\circ\phi(A)
=&\, h_k\circ\phi(\psi_k(f''_k/X_k)\circ\omega^{jk}_{(f''_j,f''_k)}) \\
=&\, h_j\circ\phi(\psi_j(f''_j/X_j)) \\
=&\, h_j\circ\phi(\psi_j(f'_j/X_j))\circ\phi(\psi_j(t/X_j)) \\
=&\, h_k\circ(\phi\circ\psi_k(f'_k/X_k))\circ\phi(B)
\end{aligned}
$$
whence the contention, since $h_k\circ(\phi\circ\psi_k(f'_k/X_k))$
is cartesian and $\phi$ is faithful. Next, for every $j,k,l\in I$
set $\cC/X_{jkl}:=\cC/X_{jk}\times_{\cC/X}\cC/X_l$; the objects of
this category are the triples
$(f'_j:X'\to X_j,f'_k:X'\to X_k,f'_l:X'\to X_l)$ such that
$f_j\circ f'_j=f_k\circ f'_k=f_l\circ f'_l$. We have obvious
projection functors
$$
\pi^0_{jkl}:\cC/X_{jkl}\to\cC/X_{kl}
\qquad
\pi^1_{jkl}:\cC/X_{jkl}\to\cC/X_{jl}
\qquad
\pi^2_{jk}:\cC/X_{jkl}\to\cC/X_{jk}.
$$
With this notation, the uniqueness properties of $\omega^{ij}$
easily imply the cocycle identity :
$$
(\omega^{jk}*\pi^2_{jkl})\odot(\omega^{kl}*\pi^0_{jkl})=
\omega^{jl}*\pi^1_{jkl}
\qquad
\text{for every $j,k,l\in I$}
$$
(details left to the reader). Let $\cS\subset\cC/X$ be the sieve
generated by the family $(f_j~|~j\in I)$; in view of proposition
\ref{prop_wecanview}, we deduce that there exist a cartesian
functor $\psi:\cS\to\cE$ and a system of isomorphisms
$(\eta_j:\psi_j\isom\psi*f_{j*}~|~j\in I)$ fulfilling the
compatibility condition
$$
(\eta_k*\pi^1_{jk})\odot\omega^{jk}=\eta_j*\pi^0_{jk}
\qquad
\text{for every $j,k\in I$}.
$$
Since $\cE$ is stack, we may then find a cartesian section
$\psi'\in\cE(X)$ with an isomorphism $\psi'_{|\cS}\isom\psi$.
Set $Y:=\psi'(\one_X)$; to conclude the proof, it will suffice
to exhibit an isomorphism $Y'\isom\phi Y$ in $\cE'$. To this
aim, pick as well a cartesian section $\tilde\psi\in\cE'(X)$
such that $\tilde\psi(\one_X)=Y'$; we set
$\tilde\psi_j:=\tilde\psi\circ f_{j*}:\cC/X_j\to\cE'$ for every
$j\in I$. Notice that
$\tilde\psi_j\circ\pi^0_{jk}=\tilde\psi_k\circ\pi^1_{jk}$ for
every $j,k\in I$. Since $h_j$ is cartesian, there exists a
unique isomorphism $h'_j:\tilde\psi(f_j)\isom\phi Y_j$ in
$F'^{-1}X_j$ such that $h_j\circ h'_j=\tilde\psi(f_j/X)$.
There follows an isomorphism of cartesian functors
$$
\tilde\eta_j:\tilde\psi_j\isom\phi\circ\psi_j
\qquad\text{such that}\qquad
\tilde\eta_{j,\one_{X_j}}=h'_j
\qquad
\text{for every $j\in I$}.
$$
Explicitly, for every $(f':X'\to X_j)\in\Ob(\cC/X_j)$, the
isomorphism
$\tilde\eta_{j,f'}:\tilde\psi_j(f')\isom\phi\circ\psi_j(f')$
is characterized by the identity :
$$
\phi(\psi_j(f'/X_j))\circ\tilde\eta_{j,f'}=
h'_j\circ\tilde\psi_j(f'/X_j)
\qquad
\text{in $\cE'$}.
$$
We claim that the following diagram commutes for every $j,k\in I$ :
$$
\xymatrix{ \tilde\psi_j\circ\pi^0_{jk} \rrdouble
\ar@{=>}[d]_{\tilde\eta_j*\pi^0_{jk}} & &
\tilde\psi_k\circ\pi^1_{jk} \ar@{=>}[d]^{\tilde\eta_k*\pi^1_{jk}} \\
\phi\circ\psi_j\circ\pi^0_{jk} \ar@{=>}[rr]^-{\phi*\omega^{ij}}
 & & \phi\circ\psi_k\circ\pi^1_{jk}.
}$$
Indeed, let $(f'_j,f'_k)$ be any object of $\cC/X_{jk}$; we compute :
$$
\begin{aligned}
h_k\circ\phi(\psi_k(f'_k/X_k))\circ
(\tilde\eta_k*\pi^1_{jk})_{(f'_j,f'_k)}
=&\, h_k\circ h'_k\circ\tilde\psi_k(f'_k/X_k) \\
=&\, \tilde\psi(f_k/X)\circ\tilde\psi(f'_k/X) \\
=&\, \tilde\psi(f_k\circ f'_k/X) \\
=&\, \tilde\psi(f_j\circ f'_j/X) \\
=&\, h_j\circ h'_j\circ\tilde\psi_j(f'_j/X_j) \\
=&\, h_j\circ\phi(\psi_j(f'_j/X_j))\circ\tilde\eta_{j,f'} \\
=&\, h_k\circ\phi(\psi_k(f'_k/X_k))\circ
\phi(\omega^{ij}_{(f'_j,f'_k)})\circ\tilde\eta_{j,f'}
\end{aligned}
$$
whence the claim, as usual, since $h_k\circ\phi(\psi_k(f'_k/X_k))$
is cartesian. In view of proposition \ref{prop_wecanview}, we
deduce that there exists a unique isomorphism of cartesian
functors
$$
\tilde\psi_{|\cS}\isom\phi\circ\psi=(\phi\circ\psi')_{|\cS}.
$$
But since $\cE'$ is $1$-separated, the latter comes from
a unique isomorphism $\tilde\psi\isom\phi\circ\psi'$, whence
the sought isomorphism $Y'\isom\phi Y$.
\end{proof}

\begin{proposition}\label{prop_unit-sep-cov-faith}
Let $(\cC,J)$ be a small site, $F:\cE\to\cC$ a fibration,
and $j\in\{0,1,2\}$. Then:
\begin{enumerate}
\item
The unit of adjunction $\eta_\cE:\cE\to\cE^a$ is $i$-covering
for $i=0,1,2$.
\item
$\cE$ is $j$-separated if and only if\/ $\eta_\cE$ is $j$-faithful.
\item
A cartesian functor $\phi:\cE\to\cE'$ of\/ $\cC$-fibrations
is $j$-covering if and only if the same holds for the induced
functor $\phi^a:\cE^a\to\cE'^a$.
\end{enumerate}
\end{proposition}
\begin{proof}(i): Quite generally, let $G:\cA\to\cB$ be a
$\cC$-cartesian functor of fibrations, and $H:\cA'\isom\cA$
a $\cC$-equivalence of categories; then it is easily seen
that $G$ is $i$-covering for some $i\in\{0,1,2\}$ if and
only if the same holds for $G\circ H$. Taking into account
lemma \ref{lem_yoga-i-coverings}(i), we are therefore
reduced to checking :

\begin{claim}\label{cl_plus-is-all-acovering}
The functor $j_\cE:\sC(\cE)\to\cE^+$ is $i$-covering for $i=0,1,2$.
\end{claim}
\begin{pfclaim} Let $[X,\psi:\cS\to\cE]$ be any object of
$\cE^+$, and pick a generating family $(f_j:X_j\to X~|~j\in I)$
for $\cS$; then $(X_j,\psi\circ f_j:\cC/X_j\to\cE)$ is an object
of $\sC(\cE)$ and $f_j$ induces a cartesian morphism
$j_\cE(X_j,\psi\circ f_{j*})\to[X,\psi]$ for every $j\in I$.
Thus, $j_\cE$ is $0$-covering.

Next, let
$[t,\sigma]:[X',\psi':\cC/X'\to\cE]\to[X,\psi:\cC/X\to\cE]$
be a morphism in $\cE^+$. By definition, $t:X'\to X$ is a
morphism of $\cC$, and
$\sigma:\psi'_{|\cS}\Rightarrow(\psi\circ t_*)_{|\cS}$ is a
natural $\cC$-transformation defined on a covering sieve
$\cS\subset\cC/X'$. Again, we pick a generating family
$(f_j:X'_i\to X'~|~j\in I)$ for $\cS$, and denote by
$(f_j,h_j):(X'_j,\psi'\circ f_{j*})\to(X',\psi')$ the
cartesian morphism in $\sC(\cE)$ induced by $f_j$. Then
$\sigma*f_{j*}:\psi'\circ f_{j*}\Rightarrow\psi\circ(t\circ f_j)_*$
defines a morphism
$(t\circ f_j,g_j):(X'_j,\psi'\circ f_{j*})\to(X,\psi)$ in
$\sC(\cE)$ such that
$j_\cE(t\circ f_j,g_j)=[t,\sigma]\circ j_\cE(f_j,h_j)$
for every $j\in I$. This proves that $j_\cE$ is $1$-covering.

Lastly, let $(t_1,\sigma_1),(t_2,\sigma_2):
(X',\psi':\cC/X'\to\cE)\to(X,\psi:\cC/X\to\cE)$ be two
morphisms in $\sC(\cE)$ such that
$j_\cE(t_1,\sigma_1)=j_\cE(t_2,\sigma_2)$. Especially, this
means that $t:=t_1=t_2$, and there exists a covering sieve
$\cS\subset\cC/X'$ such that $(\sigma_1)_{|\cS}=(\sigma_2)_{|\cS}$.
We pick a generating family $(f_j:X'_j\to X'~|~j\in I)$ for
$\cS$, and let $(f_j,h_j):(X'_j,\psi'\circ f_{j*})\to(X',\psi')$
be the cartesian morphism induced by $f_j$; then
$(t,\sigma_1)\circ(f_j,h_j)=(t,\sigma_2)\circ(f_j,h_j)$ for
every $j\in I$, so $j_\cE$ is $2$-covering.
\end{pfclaim}

(ii): From (i) and proposition \ref{prop_faith-separ-cover}
we see already that if $\cE$ is $j$-separated, then $\eta_\cE$
is $j$-faithful. Suppose then that $\eta_\cE$ is $j$-faithful;
let $X\in\Ob(\cC)$ be any object, and $\cS\subset\cC/X$ any
covering sieve; we consider the induced commutative diagram
of categories :
$$
\xymatrix{ \cE(X) \ar[r] \ar[d]_{\sCart_\cC(\cC/X,\eta_\cE)}
& \sCart_\cC(\cS,\cE) \ar[d]^{\sCart_\cC(\cS,\eta_\cE)} \\
\cE^a(X) \ar[r] & \sCart_\cC(\cS,\cE^a)
}$$
whose bottom horizontal arrow is an equivalence, since
$\cE^a$ is a stack, and whose vertical arrows are $j$-faithful,
by virtue of corollary \ref{cor_fibrations}(ii.a). It follows
that the top horizontal arrow is $j$-faithful as well, which
means that $\cE$ is $j$-separated.

(iii): We have an essentially commutative diagram of
cartesian functors
$$
\xymatrix{ \cE \ar[r]^-\phi \ar[d] & \cE' \ar[d] \\
\cE^a \ar[r]^-{\phi^a} & \cE'^a
}$$
whose vertical arrows are $i$-covering for $i=0,1,2$,
by (i). Then the assertion follows immediately from lemma
\ref{lem_yoga-i-coverings}(i,ii).
\end{proof}

\subsection{Local calculus of fractions}
\label{sec_local-calc-fract}
Let $C:=(\cC,J)$ be a site, $\cA\to\cC$ a fibration; choose
a unital cleavage $\blambda$ for $\cA$, and let
$\sc:\cC^o\to\bCat$ be the associated unital pseudo-functor.
To ease notation, for every morphism $f:X'\to X$ in $\cC$
we let $f^*:=\sc_f:\cA_X\to\cA_{X'}$, and denote by 
$\gamma^\sc_{\bullet,\bullet}$ the coherence constraint of $\sc$.
Consider a system of sets
$\Sigma_\bullet:=(\Sigma_X~|~X\in\Ob(\cC))$ with :
\begin{itemize}
\item
$\Sigma_X\subset\rMorph(\cA_X)$ for every $X\in\Ob(\cC)$
\item
$f^*(\Sigma_X)\subset\Sigma_{X'}$ for every morphism
$f:X'\to X$ of $\cC$.
\end{itemize}

\begin{proposition}\label{prop_go-for-kill}
In the situation of \eqref{sec_local-calc-fract}, there exist
a $1$-separated prestack $\cA\{\Sigma^{-1}_\bullet\}$ over $C$,
and a cartesian functor
$$
L_{\{\Sigma\}}:\cA\to\cA\{\Sigma^{-1}_\bullet\}
$$
with the following properties :
\begin{enumerate}
\item
$L_{\{\Sigma\}}f$ is an isomorphism in $\cA\{\Sigma^{-1}_\bullet\}$,
for every $f\in\Sigma:=\bigcup_{f\in\Ob(\cC)}\Sigma_X$.
\item
For every $1$-separated prestack $\cE$ on $C$, and every
cartesian functor $F:\cA\to\cE$ such that $Ff$ is invertible
in $\cE$ for every $f\in\Sigma$, there exists a unique
cartesian functor $F':\cA\{\Sigma^{-1}_\bullet\}\to\cB$
such that $F=F'\circ L_{\{\Sigma\}}$.
\end{enumerate}
\end{proposition}
\begin{proof} According to theorem \ref{th_localize-cats},
for every morphism $f:X'\to X$ of $\cC$ the functor $f^*$
extends uniquely to a functor
$f_\Sigma^*:\cA_X[\Sigma^{-1}_X]\to\cA_{X'}[\Sigma^{-1}_{X'}]$,
and $\gamma^\sc_{f,g}:f^*\circ g^*\isom(g\circ f)^*$
extends uniquely to an isomorphism of functors
$\gamma^\Sigma_{f,g}:f^*_\Sigma\circ g^*_\Sigma\isom(g\circ f)^*_\Sigma$
for every composable pair $(f,g)$ of morphisms of $\cC$ (corollary
\ref{cor_of-localization}). It is then easily seen that the
rules : $X\mapsto\cA_X[\Sigma^{-1}_X]$ for every $X\in\Ob(\cC)$
and $f\mapsto f^*_\Sigma$ for every morphism $f$ of $\cC$ yield
a unital pseudo-functor $\sc[\Sigma^{-1}_\bullet]:\cC^o\to\bCat$
with coherence constraint $\gamma^\Sigma_{\bullet\bullet}$, and we set
$$
\cA[\Sigma^{-1}_\bullet]:=\cFib(\sc[\Sigma^{-1}_\bullet]).
$$
From theorem \ref{th_localize-cats} it follows easily that
the system of localization functors $\cA_X\to\cA_X[\Sigma^{-1}_X]$
yields a unique cartesian functor $L':\cA\to\cA[\Sigma^{-1}_\bullet]$
such that the following holds. For every fibration $\cF$ over
$\cC$, and every cartesian functor $F:\cA\to\cF$ such that
$Ff$ is invertible in $\cF$ for every $f\in\Sigma$, there
exists a unique cartesian functor $F':\cA[\Sigma^{-1}_\bullet]\to\cF$
such that $F=F'\circ L'$. By theorem \ref{th_be-positive} we
conclude that the composition of $L'$ with the natural cartesian
functor $\cA[\Sigma^{-1}_\bullet]\to\cA\{\Sigma^{-1}_\bullet\}:=
\cA[\Sigma^{-1}_\bullet]^\sep$ fulfills the stated conditions.
\end{proof}

\begin{definition}\label{def_LCF}
In the situation of \eqref{sec_local-calc-fract}, we shall say
that the system $\Sigma_\bullet$
{\em admits a right local calculus of fractions} if the following
conditions hold for every $X\in\Ob(\cC)$ :
\begin{itemize}
\item[(LCF1)]
The set $\Sigma_X$ contains the isomorphisms of $\cA_X$.
\item[(LCF2)]
The set $\Sigma_X$
fulfills condition (CF2) of definition \ref{def_right-calculus}.
\item[(LCF3)]
For every morphism  $f:A\to B$ in $\cA_X$ and every
$s:C\to B$ in $\Sigma_X$, there exist a covering family
$(\phi_\lambda:X_\lambda\to X~|~\lambda\in\Lambda)$ in the
site $C$, and for every $\lambda\in\Lambda$ a morphism
$g_\lambda:D_\lambda\to\phi^*_\lambda C$ in $\cA_{X_\lambda}$ and
$t_\lambda:D_\lambda\to\phi^*_\lambda A$ in $\Sigma_{X_\lambda}$ such
that $\phi^*_\lambda(f)\circ t_\lambda=\phi^*_\lambda(s)\circ g_\lambda$.
\item[(LCF4)]
If $f,g:A\to B$ are any two morphisms in $\cA_X$ such that
$s\circ f=s\circ g$ for some $s:B\to C$ in $\Sigma_X$, then
there exist a covering family
$(\phi_\lambda:X_\lambda\to X~|~\lambda\in\Lambda)$ in the site
$C$ and for every $\lambda\in\Lambda$ an element
$t_\lambda:D_\lambda\to\phi^*_\lambda A$ in $\Sigma_{X_\lambda}$ such
that $\phi^*(f)\circ t_\lambda=\phi^*_\lambda(g)\circ t_\lambda$.
\end{itemize}
\end{definition}

\sset\subsubsection{}\label{subsec_go-for-kill}
In the situation of \eqref{sec_local-calc-fract}, for every
$X\in\Ob(\cC)$ and $A\in\Ob(\cA_X)$ let $\Sigma_{X/A}$ be the
full subcategory of $\cA_X/A$ whose objects are the elements
of $\Sigma_X$ with target equal to $A$. For every morphism
$(Y'\xrightarrow{\phi'}X)\xrightarrow{\psi/X}(Y\xrightarrow{\phi}X)$
of $\cC/X$ notice that $\psi^*:\cA_Y\to\cA_{Y'}$ induces a functor
$$
(\psi/X)^*:\cA_Y/\phi^*A\to\cA_{Y'}/\phi'^*A
\qquad
(I\xrightarrow{f}\phi^*A)\mapsto
(\psi^*I\xrightarrow{\gamma^\sc_{(\psi,\phi),A}\circ\psi^*(f)}\phi'^*A)
$$
Suppose now that $\Sigma_\bullet$ fulfills condition (LCF1) of
definition \ref{def_LCF}; then $(\psi/X)^*$ restricts to
$$
(\psi/X)^*:\Sigma_{Y/\phi^*A}\to\Sigma_{Y'/\phi'^*A}.
$$
Moreover, if $(\psi'/X):(Y''\xrightarrow{\phi''}X)\to
(Y'\xrightarrow{\phi'}X)$ is another morphism of $\cC/X$,
we have an isomorphism in $\Sigma_{Y''/\phi''^*A}$ for every
$(s:I\to\phi^*A)\in\Ob(\Sigma_{Y/\phi^*A})$
$$
(\gamma^\sc_{(\psi',\psi),I}/X):(\psi'/X)^*\circ(\psi/X)^*(s)\isom
(\psi\circ\psi'/X)^*(s).
$$
Clearly the rule $s\mapsto(\gamma^\sc_{(\psi',\psi),I}/X)$
yields an isomorphism of functors $\gamma^\Sigma_{\psi',\psi}:
(\psi'/X)^*\circ(\psi/X)^*\isom(\psi\circ\psi'/X)^*$, and
then it easily seen that the rules :
$\phi\mapsto\Sigma_{Y/\phi^*A}$ for every
$(\phi:Y\to X)\in\Ob(\cC/X)$ and $\psi/X\mapsto(\psi/X)^*$
for every morphism $\psi/X$ of $\cC/X$ define a unital
pseudo-functor $\Sigma_{\bullet/A}:(\cC/X)^o\to\bCat$, whose
coherence constraint is given by the system of isomorphisms
$\gamma^\Sigma_{\psi',\psi}$. There follows a fibration
$$
S_A:=\cFib(\Sigma_{\bullet/A})\to\cC/X
$$
and we denote by $\blambda_A$ its natural unital cleavage.
Next, for $A,B\in\Ob(\cA_X)$ consider the functor
$$
\cH_{A,B}:S_A^o\to\Set
\qquad
(Y\xrightarrow{\phi}X,I\xrightarrow{s}\phi^*A)\mapsto
\{(s,f)~|~f\in\Hom_{\cA_Y}(I,\phi^*B)\}
$$
which assigns to every morphism
$(\psi/X,h/\phi^*A):(\phi,s)\to(\phi',s')$ of $S_A$ the map
$$
\cH_{A,B}(\psi/X,h/\phi^*A):\cH_{A,B}(\phi',s')\to\cH_{A,B}(\phi,s)
\qquad
(s',f)\mapsto(s,((\psi/X)^*f)\circ h).
$$

\begin{lemma}\label{lem_S_A-loc-filtered}
With the notation of \eqref{subsec_go-for-kill}, suppose
that the system $\Sigma_\bullet$ admits a right local
calculus of fraction. Then the fibration $S_A$ is locally
cofiltered over the site $C/X$.
\end{lemma}
\begin{proof} Condition (LCF1) implies immediately condition
(a) of definition \ref{def_locally-cofiltered}. Next, let
$(\phi:Y\to X)\in\Ob(\cC/X)$, and $(\phi,s:I\to\phi^*A)$
and $(\phi,s':I'\to\phi^*A)$ be two objects of
$\Sigma_{Y/\phi^*A}$. By (LCF3) we may find a covering
$(\psi_\lambda:Y_\lambda\to Y~|~\lambda\in\Lambda)$ for the
topology $J$, and for every $\lambda\in\Lambda$ a morphism
$f_\lambda:J\to\psi_\lambda^*I$ in $\cA_{Y_\lambda}$ and an element
$t_\lambda:J\to\psi_\lambda^*I'$ of $\Sigma_{Y_\lambda}$ such that
$s''_\lambda:=\psi_\lambda^*(s)\circ f_\lambda=
\psi_\lambda^*(s')\circ t_\lambda$. By (LCF1) and (LCF2), we
have $u_\lambda:=\gamma_{(\phi,\psi_\lambda),A}\circ s''_\lambda\in
\Sigma_{Y_\lambda/(\phi\circ\psi_\lambda)^*A}$ for every
$\lambda\in\Lambda$, and it follows that $t_\lambda$ and
$f_\lambda$ yield morphisms $u_\lambda\to(\psi/X)^*(s')$
and respectively $u_\lambda\to(\psi/X)^*(s)$. Hence $S_A$
fulfills condition (b) of definition
\ref{def_locally-cofiltered}.

Lastly, consider two morphisms
$(g/\phi^*A),(g'/\phi^*):(s:I\to\phi^*A)\to(s':I\to\phi^*A)$
in $\Sigma_{Y/\phi^*A}$. In other words, $g,g':I\to I'$ are
morphisms in $\cA_Y$ with $s'\circ g=s'\circ g'=s$. By
(LCF3) there exist then a covering family
$(\psi_\lambda:Y_\lambda\to Y~|~\lambda\in\Lambda)$ and for
every $\lambda\in\Lambda$ an element
$h_\lambda:J\to\psi_\lambda^*I$ of $\Sigma_{\cA_{Y_\lambda}}$ with
$t_\lambda:=\psi^*_\lambda(g)\circ h_\lambda=
\psi^*_\lambda(g')\circ h_\lambda$. Notice that $u_\lambda:=
\psi^*_\lambda(s')\circ t_\lambda=\psi^*_\lambda(s)\circ h_\lambda$
lies in $\Sigma_{Y_\lambda,(\psi_\lambda)^*\phi^*A}$. Hence $h_\lambda$
yields a morphism $(h/(\phi\circ\psi_\lambda)^*A):
\gamma_{(\phi,\psi_\lambda),A}\circ u_\lambda\to(\psi/X)^*(s)$
in $\Sigma_{Y_\lambda,(\phi\circ\psi_\lambda)^*A}$ with
$(\psi_\lambda/X)^*(g)\circ(h/(\phi\circ\psi_\lambda)^*A)=
(\psi_\lambda/X)^*(g')\circ(h/(\phi\circ\psi_\lambda)^*A)$.
This shows that condition (c) holds as well.
\end{proof}

\sset\subsubsection{}\label{subsec_here-we-go-at-last}
With the notation of \eqref{subsec_Fubini-sheaf}, consider
the presheaf on $\cC/X$ and the sheaf on $C/X$ :
$$
M_{A,B}:=\int^{\blambda_A}\cH_{A,B}
\qquad\text{and}\qquad
H_{A,B}:=\int^{\blambda_A}_a\cH_{A,B}.
$$
Recall that for every $(\phi:Y\to X)\in\Ob(\cC/X)$, the set
$M_{AB}(\phi)$ consists of the equivalence classes of $[s,f]$
of pairs with $s\in\Ob(\Sigma_{Y/\phi^*A})$ and
$(s,f)\in\cH_{AB}(\phi,s)$. We shall denote by
$\{s,f\}\in H_{A,B}(\phi)$ the image of such a class $[s,f]$.
Taking into account lemmata \ref{lem_locally-filt-equality}
and \ref{lem_S_A-loc-filtered} we see that if $\Sigma_\bullet$
admits a right local calculus of fraction, the following holds :
\begin{itemize}
\item[(a)]
For every $\sigma\in H_{AB}(\phi)$ we have a covering
family $(\phi_\lambda:Y_\lambda\to Y~|~\lambda\in\Lambda)$ in
the site $C$, and for every $\lambda\in\Lambda$ a morphism
$g_\lambda:C_\lambda\to\phi^*_\lambda B$ in $\cA_{Y_\lambda}$ and
an element $s_\lambda:C_\lambda\to\phi^*_\lambda A$ of\/
$\Sigma_{Y_\lambda}$ such that $H_{AB}(\phi_\lambda/X)(\sigma)$
is the image of $[s_\lambda,g_\lambda]$ in
$H_{AB}(\phi\circ\phi_\lambda/X)$.
\item[(b)]
For every two pairs
$(s:I\to\phi^*A,f),(s':I'\to\phi^*A,f')\in\cH_{AB}(\phi,s)$ with
$[s,f]=[s',f']$ in $M_{AB}(\phi,s)$ there exist a covering
family $(\phi_\lambda:Y_\lambda\to Y~|~\lambda\in\Lambda)$ in
the site $C$, and for every $\lambda\in\Lambda$ two morphisms
in $\cA_{Y_\lambda}$
$$
\phi_\lambda^*I'\xleftarrow{t'_\lambda}J_\lambda
\xrightarrow{t_\lambda}\phi_\lambda^*I
\qquad\text{with}\qquad
(\phi_\lambda/X)^*(f)\circ t_\lambda=
(\phi_\lambda/X)^*(f')\circ t'_\lambda
$$
and such that $(\phi_\lambda/X)^*(s)\circ t_\lambda=
(\phi_\lambda/X)^*(s')\circ t'_\lambda$ lies in
$\Sigma_{Y_\lambda,(\phi\circ\phi_\lambda)^*A}$.
\end{itemize}
Let $\blambda_\Sigma$ and $\blambda_{\{\Sigma\}}$ be the cleavages of
$\cA[\Sigma^{-1}_\bullet]$ and respectively $\cA\{\Sigma^{-1}_\bullet\}$
deduced from the cleavage $\blambda$ of $\cA$, as in the proofs of
proposition \ref{prop_go-for-kill} and theorem \ref{th_be-positive}.
For every $A\in\Ob(\cA\{\Sigma^{-1}_\bullet\}_X)=
\Ob(\cA_X[\Sigma^{-1}_\bullet])=\Ob(\cA_X)$ define the cartesian
sections $\beta^{\blambda_\Sigma}_A:\cC/X\to\cA[\Sigma^{-1}_\bullet]$ and
$\beta^{\blambda_{\{\Sigma\}}}_{X,A}:\cC/X\to\cA\{\Sigma^{-1}_\bullet\}$
as in the proof of claim \ref{cl_split-fibrations}. Let also
$L_\Sigma:\cA\to\cA[\Sigma^{-1}_\bullet]$ be the localization
functor; we have a natural morphism of presheaves on $\cC/X$ :
$$
M_{A,B}\to\cCart(\beta^{\blambda_\Sigma}_A,\beta^{\blambda_\Sigma}_B)
\qquad
\text{for every $A,B\in\Ob(\cA_X)$}
$$
which assigns to every $[s,f]\in M_{A,B}(\phi)$ the unique
natural $\cC$-transformation
$\beta^{\blambda_\Sigma}_A\circ\phi_*\to\beta^{\blambda_\Sigma}_B\circ\phi_*$
such that $\one_Y\mapsto L_\Sigma(f)\circ L_\Sigma(s)^{-1}$, for every
$\phi\in\Ob(\cC/X)$. We deduce a morphism
$H_{A,B}\to\cCart(\beta^{\blambda_\Sigma}_A,\beta^{\blambda_\Sigma}_B)^a$
of associated sheaves on $C/X$. But recall that the proof of
theorem \ref{th_be-positive} yields as well a natural identification
$$
\cCart(\beta^{\blambda_\Sigma}_A,\beta^{\blambda_\Sigma}_B)^a\isom
\cCart(\beta^{\blambda_{\{\Sigma\}}}_A,\beta^{\blambda_{\{\Sigma\}}}_B).
$$
The composition of the latter two morphisms is a morphism
of sheaves on $C/X$ :
\set\begin{equation}\label{eq_sheaf-violante}
H_{AB}\to
\cCart(\beta^{\blambda_{\{\Sigma\}}}_{X,A},\beta^{\blambda_{\{\Sigma\}}}_{X,B})
\qquad
\{s,f\}\mapsto L_{\{\Sigma\}}(f)\circ L_{\{\Sigma\}}(s)^{-1}.
\end{equation}

\begin{proposition}
Suppose that the system $\Sigma_\bullet$ admits a right local
calculus of fraction. Then \eqref{eq_sheaf-violante} is an
isomorphism of sheaves for every $A,B\in\Ob(\cA_X)$.
\end{proposition}
\begin{proof} We define for every $X\in\Ob(\cC)$ a category
$$
\cL_X
$$
whose set of objects is $\Ob(\cA_X)$, and such that
$\Hom_{\cL_X}(A,B):=H_{AB}(\one_X)$ for every $A,B\in\Ob(\cA_X)$.
For every $A,B,C\in\Ob(\cA_X)$, the composition law
\set\begin{equation}\label{eq_boya-deh}
H_{AB}(\one_X)\times H_{BC}(\one_X)\to H_{AC}(\one_X)
\end{equation}
is obtained as follows. Let $(\phi:Y\to X)\in\Ob(\cC/X)$
and $[s:I\to\phi^*A,f]\in M_{AB}(\phi)$ and
$[s':I'\to\phi^*B,f']\in M_{BC}(\phi)$ any two sections;
by virtue of (LCF3) we may find a covering family
$(\phi_\lambda:Y_\lambda\to Y~|~\lambda\in\Lambda)$
and for every $\lambda\in\Lambda$, morphisms in $\cA_{Y_\lambda}$ :
$$
\phi^*_\lambda I\xleftarrow{s_\lambda}J_\lambda
\xrightarrow{f_\lambda}\phi^*_\lambda I'
\qquad\text{such that}\qquad
(\phi_\lambda/X)^*(f)\circ s_\lambda=
(\phi_\lambda/X)^*(s')\circ f_\lambda
\quad\text{and}\quad
s_\lambda\in\Sigma_{Y_\lambda}.
$$
Set $t_\lambda:=(\phi_\lambda/X)^*(s)\circ s_\lambda:
J_\lambda\to(\phi\circ\phi_\lambda)^*A$ and
$g_\lambda:=(\phi_\lambda/X)^*(f')\circ f_\lambda:
J_\lambda\to(\phi\circ\phi_\lambda)^*B$ for every $\lambda$.
For $\lambda,\mu\in\Lambda$, consider morphisms
$Y_\lambda\xleftarrow{\psi_\lambda}Z\xrightarrow{\psi_\mu}Y_\mu$
such that $\psi:=\phi_\lambda\circ\psi_\lambda=\phi_\mu\circ\psi_\mu$.

\begin{claim}\label{cl_boia-deh}
$H_{AB}(\psi_\lambda/X)\{t_\lambda,g_\lambda\}=
H_{AB}(\psi_\mu/X)\{t_\mu,g_\mu\}$.
\end{claim}
\begin{pfclaim} We need to show that
$\{(\psi_\lambda/\!X)^*t_\lambda,(\psi_\lambda/\!X)^*g_\lambda\}=
\{(\psi_\mu/\!X)^*t_\mu,(\psi_\mu/\!X)^*g_\mu\}$. However, by (LCF3)
we may find a covering family
$(\rho_\nu:Z_\nu\to Z~|~\nu\in\Lambda')$ and for every
$\nu\in\Lambda'$, morphisms $\rho_\nu^*\psi_\mu^*J_\mu
\xleftarrow{s_\nu}J_\nu\xrightarrow{s'_\nu}
\rho_\nu^*\psi_\lambda^*J_\lambda$ in $\cA_{Z_\nu}$ with
$s_\nu\in\Sigma_{Z_\nu}$, and such that
$$
(\rho_\nu/Y)^*(\psi_\mu/Y)^*(s_\mu)\circ s_\nu=
(\rho_\nu/Y)^*(\psi_\lambda/Y)^*(s_\lambda)\circ s'_\nu.
$$
We compute :
$$
\begin{aligned}
(\rho_\nu/Y)^*((\psi/X)^*(s')\circ(\psi_\mu/Y)^*(f_\mu))\circ s_\nu&=
(\rho_\nu/Y)^*((\psi_\mu/X)^*((\phi_\mu/X)^*(s')\circ f_\mu))\circ s_\nu \\
&=(\rho_\nu/Y)^*((\psi_\mu/X)^*((\phi_\mu/X)^*(f)\circ s_\lambda))\circ s_\nu \\
&=(\rho_\nu/Y)^*((\psi/X)^*(f)\circ(\psi_\mu/Y)^*(s_\mu))\circ s_\nu \\
&=(\rho_\nu/Y)^*((\psi/X)^*(f)\circ(\psi_\lambda/Y)^*(s_\lambda))\circ s'_\nu \\
&=(\rho_\nu/Y)^*((\psi/X)^*(s')\circ(\psi_\lambda/Y)^*(f_\lambda))\circ s'_\nu.
\end{aligned}
$$
By (LCF4) it follows that for every $\nu\in\Lambda'$ there
exist a covering family
$(\rho_{\nu\nu'}:Z_{\nu\nu'}\to Z_\nu~|~\nu'\in\Lambda'_\nu)$
and for every $\nu'\in\Lambda'_\nu$ an element
$t_{\nu\nu'}:J_{\nu\nu'}\to\rho_{\nu\nu}^*J_\nu$ of
$\Sigma_{Z_{\nu\nu'}}$ such that
$$
(\rho_{\nu\nu'}/Y)^*(\rho_\nu/Y)^*
((\psi_\mu/Y)^*(f_\mu)\circ s_\nu)\circ t_{\nu\nu'}=
(\rho_{\nu\nu'}/Y)^*(\rho_\nu/Y)^*
((\psi_\lambda/Y)^*(f_\lambda)\circ s'_\nu)
\circ t_{\nu\nu'}.
$$
Set $\rho'_{\nu\nu'}:=\rho_\nu\circ\rho_{\nu\nu'}$ for every
$\nu'\in\Lambda'_\nu$, and
$$
T_{\lambda\nu\nu'}:=(\psi_\lambda\circ\rho'_{\nu\nu'}/X)^*t_\lambda
\qquad
G_{\lambda\nu\nu'}:=(\psi_\lambda\circ\rho'_{\nu\nu'}/X)^*g_\lambda
\qquad
\gamma_{\lambda\nu\nu'}:=\gamma^\sc_{(\rho'_{\nu\nu'},\psi_\lambda),J_\lambda}
$$
and define likewise $T_{\mu\nu\nu'}$, $G_{\mu\nu\nu'}$ and
$\gamma_{\mu\nu\nu'}$.
We compute :
$$
\begin{aligned}
[T_{\lambda\nu\nu'},G_{\lambda\nu\nu'}]&=
[T_{\lambda\nu\nu'}\circ\gamma_{\lambda\nu\nu'}\circ
(\rho_{\nu\nu'}/Z)^*(s'_\nu)\circ t_{\nu\nu'},
G_{\lambda\nu\nu'}\circ\gamma_{\lambda\nu\nu'}\circ
(\rho'_{\nu\nu'}/Z)^*(s'_\nu)\circ t_{\nu\nu'}] \\
&=[T_{\mu\nu\nu'}\circ\gamma_{\mu\nu\nu'}\circ
(\rho_{\nu\nu'}/Z)^*(s_\nu)\circ t_{\nu\nu'},
G_{\mu\nu\nu'}\circ\gamma_{\mu\nu\nu'}\circ
(\rho_{\nu\nu'}/Z)^*(s_\nu)\circ t_{\nu\nu'}] \\
&=[T_{\mu\nu\nu'},G_{\mu\nu\nu'}]
\end{aligned}
$$
whence the contention.
\end{pfclaim}

Claim \ref{cl_boia-deh} implies that the system
$(\{t_\lambda,g_\lambda\}~|~\lambda\in\Lambda)$ defines
a unique section of $H_{AC}(\phi)$ that we denote
$\{s',f'\}\circ\{s,f\}$. Notice that the argument proves
as well that the construction of $\{s',f'\}\circ\{s,f\}$
is independent of all auxiliary choices. A direct inspection
shows then that we obtain a well defined morphism of presheaves
$$
M_{AB}\times M_{BC}\to H_{AC}
\qquad
([s,f],[s',f'])\mapsto\{s',f'\}\circ\{s,f\}
$$
which therefore induces a unique morphism of sheaves on
$C/X$ :
$$
H_{AB}\times H_{BC}\to H_{AC}
\qquad\text{such that}\qquad
(\{s,f\},\{s',f'\})\mapsto\{s',f'\}\circ\{s,f\}.
$$
Then the sought map \eqref{eq_boya-deh} is the map of
sets obtained by evaluating the foregoing morphism on
the object $\one_X\in\Ob(\cC/X)$. Let us check that
$\{\one_{\phi^*B},\one_{\phi^*B}\}\circ\sigma=\sigma=
\sigma\circ\{\one_{\phi^*A},\one_{\phi^*A}\}$ for every
$A,B\in\Ob(\cA_X)$ and every $\phi\in\Ob(\cC/X)$ :
indeed, it suffices to check these identities for
$\sigma=\{s,f\}$, for any $(s,f)\in\cH_{A,B}(\phi)$,
and the latter follow by simple inspection. Lastly,
in order to show the associativity property for our
composition law, consider $A,B,C,D\in\Ob(\cA_X)$
and $(s:I\to\phi^*A,f),(s':I'\to\phi^*B,f'),
(s'':I''\to\phi^*C,f'')\in\cH_{AB}(\phi)$, for some
$(\phi:Y\to X)\in\Ob(\cC/X)$. By (LCF3) we may find
a covering family
$(\phi_\lambda:Y_\lambda\to Y~|~\lambda\in\Lambda)$ and
for every $\lambda\in\Lambda$, morphisms in $\cA_{Y_\lambda}$ :
$$
\phi^*_\lambda I\xleftarrow{t_\lambda}J_\lambda
\xrightarrow{g_\lambda}\phi^*_\lambda I'
\xleftarrow{t'_\lambda}J'_\lambda
\xrightarrow{g'_\lambda}\phi^*_\lambda I''
\qquad\text{with}\qquad
t_\lambda,t'_\lambda\in\Sigma_{Y_\lambda}
$$
such that $(\phi_\lambda/X)^*(f)\circ t_\lambda=
(\phi_\lambda/X)^*s'\circ g_\lambda$ and
$(\phi_\lambda/X)^*(f')\circ t'_\lambda=
(\phi_\lambda/X)^*(s'')\circ g'_\lambda$. Then we may also find
for every $\lambda\in\Lambda$ a covering family
$(\phi_{\lambda\mu}:Y_{\lambda\mu}\to Y_\lambda~|~\mu\in\Lambda_\lambda)$
and for every $\mu\in\Lambda_\lambda$, morphisms in
$\cA_{Y_{\lambda\mu}}$ :
$$
\phi^*_{\lambda\mu}J_\lambda\xleftarrow{t_{\lambda\mu}}
J_{\lambda\mu}\xrightarrow{g_{\lambda\mu}}\phi^*_{\lambda\mu}J'_\lambda
\qquad\text{such that}\qquad
(\phi_{\lambda\mu}/Y)^*(g_\lambda)\circ t_{\lambda\mu}=
(\phi_{\lambda\mu}/Y)^*(t_\lambda)\circ g_{\lambda\mu}
$$
and with $t_{\lambda\mu}\in\Sigma_{Y_{\lambda\mu}}$. For every
such $\lambda$ and $\mu$ set
$\psi_{\lambda\mu}:=\phi_\lambda\circ\phi_{\lambda\mu}$ and
$$
(t'_{\lambda\mu},g'_{\lambda\mu}):=
((\psi_{\lambda\mu}/X)^*(s)\circ
(\phi_{\lambda\mu}/Y)^*(t_\lambda)\circ t_{\lambda\mu},
(\psi_{\lambda\mu}/X)^*(f'')\circ
(\phi_{\lambda\mu}/Y)^*(g_\lambda)\circ g_{\lambda\mu}).
$$
It then follows easily that the system
$(\{t'_{\lambda\mu},g'_{\lambda\mu}\}~|~
\lambda\in\Lambda,\mu\in\Lambda_\lambda)$ represents both
$\{s'',f''\}\circ(\{s',f'\}\circ\{s,f\})$ and
$(\{s'',f''\}\circ\{s',f'\})\circ\{s,f\}$. The sought
associativity property is a straightforward consequence.
This concludes the construction of $\cL_X$.

Next, notice that for every morphism $\phi:Y\to X$
in $\cC$ and every $A,B\in\Ob(\cA_X)$ we have a natural
isomorphism of presheaves on $\cC/Y$ :
$$
M_{\phi^*A,\phi^*B}\isom(\phi_*)^\wedge M_{A,B}
$$
Namely, for every $(\psi:Z\to Y)\in\Ob(\cC/Y)$, the
isomorphism is given by the map :
$$
M_{\phi^*A,\phi^*B}(\psi)\isom M_{A,B}(\phi\circ\psi)
\qquad
[s,f]\mapsto
[\gamma^\sc_{(\psi,\phi),A}\circ s,\gamma^\sc_{(\psi,\phi),B}\circ f]
$$
Since $\phi_*$ is continuous and cocontinuous for the sites
$C/X$ and $C/Y$ (remark \ref{rem_continue-local}(i)), combining
with lemma \ref{lem_improve}(ii) there follows a natural
isomorphism of sheaves on $C/Y$ :
\set\begin{equation}\label{eq_keep-cool}
H_{\phi^*A,\phi^*B}\isom j^*_\phi H_{A,B}.
\end{equation}
We deduce a map
$$
\sd_{\phi,(A,B)}:\Hom_{\cL_X}(A,B)\to H_{A,B}(\phi)\isom
\Hom_{\cL_Y}(\phi^*A,\phi^*B)
$$
namely the composition of the restriction map
$H_{A,B}(\phi/X):H_{A,B}(\one_X)\to H_{A,B}(\phi)$ with the
evaluation at $\one_Y\in\Ob(\cC/Y)$ of the inverse of the
foregoing isomorphism of sheaves. A simple inspection
shows that the rules : $A\mapsto\phi^*A$ and
$\sigma\mapsto\sd_{\phi,(A,B)}(\sigma)$ for every
$A,B\in\cA_X$ and every $\sigma\in\Hom_{\cL_X}(A,B)$
yield a well defined functor
$$
\sd_\phi:\cL_X\to\cL_Y.
$$
Furthermore, for every pair
$Z\xrightarrow{\psi}Y\xrightarrow{\phi}X$ of morphisms
of $\cC$, the natural isomorphism
$\gamma^\sc_{\psi,\phi}:\psi^*\circ\phi^*\isom(\phi\circ\psi)^*$
induces a natural isomorphism
$$
\gamma^\sd_{\psi,\phi}:
\sd_\psi\circ\sd_\phi\isom\sd_{\phi\circ\psi}
\qquad
A\mapsto\{\one_{\psi^*\phi^*A},\gamma^\sc_{(\psi,\phi),A}\}.
$$
Clearly the rules : $X\mapsto\cL_X$ and $\phi\mapsto\sd_\phi$
for every $X\in\Ob(\cC)$ and every morphism $\phi$ of $\cC$
define a unital pseudo-functor $\sd:\cC^o\to\bCat$, with
coherence constraint $\gamma^\sd_{\bullet,\bullet}$. We set :
$$
\cL:=\cFib(\sd)
$$
and we denote by $\blambda_\cL$ the natural cleavage of $\cL$.
Then, for every $X\in\Ob(\cC)$ and every $A\in\Ob(\cL_X)$
let $\beta^{\blambda_\cL}_{X,A}:\cC/X\to\cL$ be the cartesian
section defined as in the proof of claim
\ref{cl_split-fibrations}; by inspecting the constructions,
we see that for every $A,B\in\Ob(\cL_X)$, the system of
isomorphisms \eqref{eq_keep-cool} induces a natural isomorphism
of presheaves on $\cC/X$ :
$$
H_{AB}\isom\cCart(\beta^{\blambda_\cL}_{X,A},\beta^{\blambda_\cL}_{X,B}).
$$
Especially, $\cCart(\sigma,\sigma')$ is a sheaf on $C/X$ for
every pair of cartesian sections $\sigma,\sigma'\in\cL(X)$,
and therefore $\cL$ is a $1$-separated prestack (lemma
\ref{lem_crit-0-1_separation}). We have a natural functor
$$
F_X:\cA_X\to\cL_X
\qquad
\text{for every $X\in\Ob(\cC)$}
$$
that is the identity on objects and is given by the rule :
$(f:A\to B)\mapsto\{\one_A,f\}$ on morphisms $f$ of $\cA_X$.
Clearly the rule $X\mapsto F_X$ defines a strict pseudo-natural
transformation $F_\bullet:\sc\Rightarrow\sd$, whence a cartesian
functor
$$
L:=\cFib(F_\bullet):\cA\to\cL.
$$
By construction, $L(s)$ is an isomorphism in $\cL$, for every
$s\in\bigcup_{X\in\Ob(\cC)}\Sigma_X$, hence $L$ factors uniquely
through a cartesian functor $L':\cA\{\Sigma^{-1}_\bullet\}\to\cL$
and the universal cartesian functor $L_{\{\Sigma\}}$ of proposition
\ref{prop_go-for-kill}. To construct conversely a functor
$\cL\to\cA\{\Sigma^{-1}_\bullet\}$, we begin by evaluating at
$\one_X$ the morphism of sheaves \eqref{eq_sheaf-violante},
to get a map
$$
L''_{X,A,B}:\Hom_{\cL_X}(A,B)\to\Hom_{\cA\{\Sigma^{-1}_\bullet\}_X}(A,B).
$$
Then, a direct inspection of the constructions yields a 
functor $L''_X:\cL_X\to\cA\{\Sigma^{-1}_\bullet\}_X$ that
is the identity on objects, and is given on morphisms by the
rule : $\sigma\mapsto L''_{X,A,B}(\sigma)$ for every
$A,B\in\Ob(\cL_X)$, and every $\sigma\in\Hom_{\cL_X}(A,B)$ : the
details shall be left to the reader.

Lastly, if $\sc_{\{\Sigma\}}$ denotes the pseudo-functor
associated with the cleavage $\blambda_{\{\Sigma\}}$, we
easily see that the rule : $X\mapsto\cL''_X$ for every
$X\in\Ob(\cC)$ defines a strict pseudo-natural
$\cC$-transformation $L''_\bullet:\sd\Rightarrow\sc_{\{\Sigma\}}$,
whence a cartesian functor
$$
L'':=\cFib(L''_\bullet):\cL\to\cA\{\Sigma^{-1}_\bullet\}.
$$
A direct inspection shows that $L'\circ L''=\one_\cL$. In
order to show that $L''\circ L'=\one_{\cA\{\Sigma^{-1}_\bullet\}}$,
it suffices to check that
$L''\circ L'\circ L_{\{\Sigma\}}=L_{\{\Sigma\}}$, due to the
universal property of $L_{\{\Sigma\}}$. Thus, we come down
to the proving that $L''\circ L=L_{\{\Sigma\}}$, which again
holds by direct inspection. Summing up, we have proven that
$L''$ is an isomorphism of categories; the proposition follows
at once.
\end{proof}

\subsection{Functorial properties of the categories of stacks}
\label{sec_funct-stacks}
In section \ref{sec_cont-functors} we have attached to
every continuous or cocontinuous functor between sites
certain natural functors on the corresponding categories
of sheaves. Hereafter we shall likewise associate with
such functors certain natural pseudo-functors on the
corresponding $2$-categories of stacks.

\begin{proposition}\label{prop_cocont-and-stacks}
Let $(\cC,J)$ and $(\cC',J')$ be two sites, $u:\cC\to\cC'$
a cocontinuous functor, and $i\in\{0,1,2\}$. The following holds :
\begin{enumerate}
\item
For every $i$-covering cartesian functor $\cE_1\xrightarrow{\phi}\cE_2$
of\/ $\cC'$-fibrations $\cE_1\xrightarrow{F_1}\cC'\xleftarrow{F_2}\cE_2$,
the functor $\Fib(u)^*(\phi):\cC\times_{\cC'}\cE_1\to\cC\times_{\cC'}\cE_2$
is $i$-covering (notation of remark {\em\ref{rem_added-little-extra}(i)}).
\item
For every $i$-separated fibration $\cF\to\cC$, the
fibration $\Fib(u)_*(\cF)$ is $i$-separated.
\end{enumerate}
\end{proposition}
\begin{proof}(i): Recall that the objects of $\cC\times_{\cC'}\cE_j$
for $j=1,2$ are the pairs $(X,Y)$ with $X\in\Ob(\cC)$,
$Y\in\Ob(\cE_j)$, such that $uX=F_jY$. The morphisms
$(X,Y)\to(X',Y')$ are the pairs $(f,g)$ where $f:X\to X'$
is a morphism in $\cC$, $g:Y\to Y'$ is a morphism in $\cE_j$,
and $u(f)=F_j(g)$.

Suppose first that $i=0$, and let $(X,Y)\in\Ob(\cC\times_{\cC'}\cE_2)$
be any object. By assumption, there exist a covering family
$(f'_j:X'_j\to uX~|~j\in I)$ and for every $j\in I$ a cartesian
morphism $h_j:\phi Y_j\to Y$ in $\cE_2$ with $F_2h_j=f'_j$. Let
$\cS\subset\cC'/uX$ be the sieve generated by $(f'_j~|~j\in I)$;
since $u$ is cocontinuous, there exists a covering family
$(f_k:X_k\to X~|~k\in I')$ such that $u(f_k)\in\Ob(\cS)$
for every $k\in I'$. Then, for every $k\in I'$ pick
$j(k)\in I$ such that $u(f_k)=f'_{j(k)}\circ f''_k$ for some
morphism $f''_k:uX_k\to X'_{j(k)}$ in $\cC'$, and choose a
cartesian morphism $h'_k:Z_k\to Y_{j(k)}$ in $\cE_1$ such that
$F_1h'_k=f''_k$. Then
$h''_k:=h_{j(k)}\circ\phi(h'_k):\phi Z_k\to Y$ is cartesian
and $F_2h''_k=u(f_k)$; therefore
$(f_k,h''_k):(X_k,\phi Z_k)\to(X,Y)$ is a cartesian morphism
in $\cC\times_{\cC'}\cE_2$ for every $j\in I$. This proves that
$\Fib(u)^*(\phi)$ is $0$-covering.

Next, suppose that $i=1$, and let
$(\one_X,g):(X,\phi Y_1)\to(X,\phi Y_2)$ be a morphism in
$\cC\times_{\cC'}\cE_2$. By assumption there exist a covering
family $(f'_j:X'_j\to uX~|~j\in I)$ in $\cC'$ and morphisms
$Y_2\xleftarrow{g_j}Z_j\xrightarrow{h_j}Y_1$ in $\cE_1$ with
$h_j$ cartesian, such that $\phi(g_j)=g\circ\phi(h_j)$ and
$F_1h_j=f'_j$ for every $j\in I$. Then, choose a covering
family $(f_k:X_k\to X~|~k\in I')$ as in the foregoing, so
that for every $k\in I'$ there exist $j(k)\in I$ and a
factorization $u(f_k)=f'_{j(k)}\circ f''_k$ for some morphism
$f''_k:uX_k\to X'_{j(k)}$ in $\cC'$. Choose also for every
$k\in I'$ a cartesian morphism $h'_k:Z'_k\to Z_{j(k)}$ with
$F_1h'_k=f''_k$. Then $h''_k:=h_{j(k)}\circ h'_k:Z'_k\to Y_1$
is cartesian and $F_1h''_k=u(f_k)$, so we get a cartesian
morphism $(f_k,h''_k):(X_k,Z'_k)\to(X,Y_1)$ such that
$(\one_X,g)\circ(f_k,\phi(h''_k))=
(f_k,\phi(g_{j(k)}\circ h'_k))$ for every $k\in I'$.
This shows that $\Fib(u)^*(\phi)$ is $1$-covering.

Lastly, say that $i=2$, and let
$(\one_X,g_1),(\one_X,g_2):(X,Y_1)\to(X,Y_2)$ be two morphisms
in $\cC\times_{\cC'}\cE_1$ with $\phi(g_1)=\phi(g_2)$. By
assumption there exist a covering family
$(f'_j:X'_j\to uX~|~j\in I)$ and a cartesian morphism
$h_j:Z_j\to Y_1$ in $\cE_1$ such that
$g_1\circ h_j=g_2\circ h_j$ and $F_1h_j=f'_j$ for every
$j\in I$. We pick again a covering family
$(f_k:X_k\to X~|~k\in I')$ in $\cC$ such that for every
$k\in I'$ there exist $j(k)\in I$ and a factorization
$u(f_k)=f'_{j(k)}\circ f''_k$, and we choose for every
$k\in I'$ a cartesian morphism $h'_k:Z'_k\to Z_{j(k)}$
such that $F_1h'_k=f''_k$. Then
$h''_k:=h_{j(k)}\circ h'_k:Z'_k\to Y_1$ is cartesian and
$F_1h''_k=u(f_k)$; therefore
$(f_k,h''_k):(X_k,Z'_k)\to(X,Y_1)$ is cartesian morphism
in $\cC\times_{\cC'}\cE_1$ such that
$(\one_X,g_1)\circ(f_k,h''_k)=(\one_X,g_2)\circ(f_k,h''_k)$
for every $k\in I'$. This proves that $\Fib(u)^*(\phi)$
is $2$-covering.

(ii): See \eqref{subsec_adjoints-of-Fib(rho)} for the definition
of the pseudo-functor $\Fib(u)_*$. Let $X\in\Ob(\cC)$ be any object,
and $\cS\subset\cC/X$ a covering sieve; we consider the essentially
commutative diagram :
$$
\xymatrix@R-5pt{ \Fib(u)_*\cF(X) \ar[r] \ar[d] &
\sCart_\cC(\Fib(u)^*(\cC/X),\cF) \ar[r] \ar[d] &
\sCart_\cC((\Fib(u)^*(\cC/X))^a,\cF^a) \ar[d] \\
\sCart_\cC(\cS,\Fib(u)_*\cF) \ar[r] & \sCart_\cC(\Fib(u)^*(\cS),\cF)
\ar[r] & \sCart_\cC((\Fib(u)^*(\cS))^a,\cF^a).
}$$
We need to show that the left vertical arrow is $i$-faithful,
and we know that the horizontal arrows of the left square
subdiagram are equivalences. We are thus reduced to checking
that the central vertical arrow is $i$-faithful. But since
$\cF$ is $i$-separated, the unit of adjunction $\cF\to\cF^a$
is $i$-faithful (proposition \ref{prop_unit-sep-cov-faith}(ii)),
so the same holds for the horizontal arrows of the right square
subdiagram (corollary \ref{cor_fibrations}(ii.a)), and we are
further reduced to checking that the right vertical arrow is
an equivalence. To this aim, it suffices to check that the
inclusion functor $t:\cS\to\cC/X$ induces an equivalence
$(\Fib(u)^*(\cS))^a\isom(\Fib(u)^*(\cC/X))^a$ (corollary
\ref{cor_fibrations}(ii.b)). Now, recall that $\cC/X$ is
isomorphic to $\cFib(h_X)$, where $h_X$ is the presheaf
on $\cC$ represented by $X$; likewise, $\cS$ is isomorphic
to $\cFib(h_\cS)$, for a covering subobject $s:h_\cS\to h_X$.
We deduce natural isomorphisms :
$$
(\Fib(u)^*(\cS))^a\isom\cFib(u^\wedge h_\cS)^a\isom
\cFib((u^\wedge h_\cS)^a)
$$
and likewise $(\Fib(u)^*(\cC/X))^a\isom\cFib((u^\wedge h_X)^a)$
(theorem \ref{th_stackfication} and remark \ref{rem_mathsfev}(ii)).
Clearly, these isomorphisms identify $s$ with $\cFib((u^\wedge s)^a)$;
the latter is an isomorphism, according to lemma
\ref{lem_breve}, whence the contention.
\end{proof}

\begin{definition}\label{def_weak-morph-of-sites}
Let $C:=(\cC,J)$ and $C':=(\cC',J')$ be two sites, and
$u:\cC\to\cC'$ a functor. We say that $u:C'\to C$ is a
{\em weak morphism of sites\/} if for every universe
$\sV$ and every $\sV$-stack $\cE$ on $C'$, the fibration
$\sV\tdu\Fib(u)^*\cE$ is a $\sV$-stack on $C$. Then, for
every such $\sV$, the pseudo-functor $\sV\tdu\Fib(u)^*$
induces by restriction a (strict) pseudo-functor
$$
\sV\tdu\St(u)_*:\sV\tdu\Stack(C')\to\sV\tdu\Stack(C).
$$
\end{definition}

\begin{remark}
(i)\ \
From example \ref{ex_sheaves-as-stacks}, it follows that
every weak morphism of sites $C'\to C$ is a continuous
functor for the topologies of $C$ and $C'$.

(ii)\ \
We shall see hereafter that every morphism of sites is
a weak morphism of sites.
\end{remark}

\begin{proposition}\label{prop_weak-morph-of-sites}
Let $C:=(\cC,J)$ and $C':=(\cC',J')$ be two sites,
$u:\cC\to\cC'$ a functor, and $\sV$ a universe such that
$\cC$ is $\sV$-small and $\cC'$ has $\sV$-small $\Hom$-sets.
The following conditions are equivalent :
\begin{enumerate}
\alphaenu
\item
For every $\sV$-stack $\cE$ on $C'$, the fibration
$\sV\tdu\Fib(u)^*\cE$ is a $\sV$-stack on $C$.
\item
For every morphism $\phi:\cE\to\cF$ in $\sV\tdu\Fib(\cC)$
that is $i$-covering for $i=0,1,2$, the morphism
$\sV\tdu\Fib(u)_!\phi:
\sV\tdu\Fib(u)_!\cE\to\sV\tdu\Fib(u)_!\cF$ is $i$-covering
for $i=0,1,2$.
\item
For every $X\in\Ob(\cC)$ and every covering sieve
$\cS\subset\cC/X$, the induced morphism of $\cC'$-fibrations
$\sV\tdu\Fib(u)_!\cS\to\sV\tdu\Fib(u)_!(\cC/X)$ is $i$-covering
for $i=0,1,2$.
\item
$u:C'\to C$ is a weak morphism of sites.
\end{enumerate}
\end{proposition}
\begin{proof}(a)$\Rightarrow$(b): In view of remark
\ref{rem_explicit-Fib_*}(v), we may replace $\sV$ by
a larger universe, and assume that both $\cC$ and $\cC'$
are $\sV$-small. Then, according to propositions
\ref{prop_faith-separ-cover} and \ref{prop_unit-sep-cov-faith}(iii),
the functor $\phi^a:\cE^a\to\cF^a$ is an equivalence, and we need
to check that $(\sV\tdu\Fib(u)_!\phi)^a$ is an equivalence.
To this aim, for every $\sV$-stack $\cA$ on $C'$, we consider
the diagram :
$$
\xymatrix@R+5pt{\sCart_{\cC'}(\Fib(u)_!(\cF)^a,\cA) \ar[r]
\ar[d]|{\sCart_{\cC'}(\Fib(u)_!(\phi)^a,\cA)} &
\sCart_{\cC'}(\Fib(u)_!(\cF),\cA) \ar[r]
\ar[d]|{\sCart_{\cC'}(\Fib(u)_!(\phi),\cA)} &
\sCart_\cC(\cF,\Fib(u)^*\cA) \ar[d]|{\sCart_\cC(\phi,\Fib(u)^*\cA)} \\
\sCart_{\cC'}(\Fib(u)_!(\cE)^a,\cA) \ar[r] &
\sCart_{\cC'}(\Fib(u)_!(\cE),\cA) \ar[r] & \sCart_\cC(\cE,\Fib(u)^*\cA)
}$$
whose left square subdiagram is induced by the natural
cartesian functors
$$
\Fib(u)_!(\cF)\to\Fib(u)_!(\cF)^a
\qquad\text{and}\qquad
\Fib(u)_!(\cE)\to\Fib(u)_!(\cE)^a
$$
and whose right square subdiagram is given by the coherence
constraint of the $2$-adjunction for the pair
$(\sV\tdu\Fib(u)_!,\sV\tdu\Fib(u)^*)$. Thus, both subdiagrams
are essentially commutative, and the same then holds for their
composition. Moreover, all horizontal arrows are equivalences,
so the left-most vertical arrow is an equivalence if and only
if the same holds for the right-most one, and in light of lemma
\ref{lem_long-time}(iii) we are reduced to checking that
$\sCart_{\cC'}(\phi,\Fib(u)^*\cA)$ is an equivalence for every
such $\cA$. But by assumption, $\sV\tdu\Fib(u)^*\cA$ is a stack
on $C$, so we have as well the pseudo-commutative diagram
$$
\xymatrix{
\sCart_\cC(\cF^a,\Fib(u)^*\cA) \ar[d]_{\sCart_\cC(\phi^a,\Fib(u)^*\cA)}
\ar[r] & \sCart_\cC(\cF,\Fib(u)^*\cA) \ar[d]^{\sCart_\cC(\phi,\Fib(u)^*\cA)}
\\
\sCart_\cC(\cE^a,\Fib(u)^*\cA) \ar[r] &
\sCart_\cC(\cE,\Fib(u)^*\cA)
}$$
whose horizontal arrows are equivalences. Hence it suffices to
show that $\sCart_\cC(\phi^a,\Fib(u)^*\cA)$ is an equivalence
for every $\sV$-stack $\cA$ on $C'$, which follows by invoking
again lemma \ref{lem_long-time}(iii).

(b)$\Rightarrow$(c): This is clear from propositions
\ref{prop_faith-separ-cover} and \ref{prop_unit-sep-cov-faith}(iii),
after inspecting the pseudo-commutative diagram \ref{eq_free-palestine}.

(c)$\Rightarrow$(d): Let $\sV'$ be any universe; we need to
check that $\sV'\tdu\Fib(u)^*\cE$ is a $\sV'$-stack on $C$,
for every $\sV'$-stack $\cE$ on $C'$. To this aim, we can
replace $\sV'$ by a larger universe, and assume that
$\sV\subset\sV'$, and both $\cC$ and $\cC'$ are $\sV'$-small.
Moreover, in light of remark \ref{rem_explicit-Fib_*}(v) and
our assumption (c), we see that for every $X\in\Ob(\cC)$ and
every covering sieve $\cS\subset\cC/X$, the induced morphism
of $\cC'$-fibrations
$\sV'\tdu\Fib(u)_!\cS\to\sV'\tdu\Fib(u)_!(\cC/X)$ is $i$-covering
for $i=0,1,2$. Now, the coherence constraint for the $2$-adjoint
pair $(\sV'\tdu\Fib(u)_!,\sV'\tdu\Fib(u)^*)$ yields an essentially
commutative diagram :
$$
\xymatrix{ \cCart_\cC(\cC/X,\sV'\tdu\Fib(u)^*\cE) \ar[r] \ar[d] &
\cCart_\cC(\cS,\sV'\tdu\Fib(u)^*\cE) \ar[d] \\
\cCart_{\cC'}(\sV'\tdu\Fib(u)_!(\cC/X),\cE) \ar[r] &
\cCart_{\cC'}(\sV'\tdu\Fib(u)_!\cS,\cE)
}$$
whose vertical arrow are equivalences. Thus, we are reduced to
checking that the bottom horizontal arrow is an equivalence.
But since $\cE$ is a stack, we have as well the following
pseudo-commutative commutative diagram, whose vertical arrows
are again equivalences :
$$
\xymatrix{ \cCart_{\cC'}(\sV'\tdu\Fib(u)_!(\cC/X)^a,\cE) \ar[r]
\ar[d] & \cCart_{\cC'}(\sV'\tdu\Fib(u)_!(\cS)^a,\cE) \ar[d] \\
\cCart_{\cC'}(\sV'\tdu\Fib(u)_!(\cC/X),\cE) \ar[r] &
\cCart_{\cC'}(\sV'\tdu\Fib(u)_!\cS,\cE)
}$$
Hence, it suffices to check that the top horizontal arrow
of this latter diagram is an equivalence. But the functor
$\sV'\tdu\Fib(u)_!(\cS)^a\to\sV'\tdu\Fib(u)_!(\cC/X)^a$ is
an equivalence, by propositions
\ref{prop_faith-separ-cover} and \ref{prop_unit-sep-cov-faith}(iii);
the contention is an immediate consequence.

(d)$\Rightarrow$(a) is trivial.
\end{proof}

\sset\subsubsection{}\label{subsec_weak-strong-continuity}
Let $C:=(\cC,J)$ and $C':=(\cC',J')$ be two sites; we wish to
give some useful criteria for a functor $u:\cC\to\cC'$ to be
a weak morphism of sites $C'\to C$. To this aim, we consider
the following conditions :
\begin{itemize}
\item[(C0)]
For every $X\in\Ob(\cC)$ and every covering family
$(X_i\to X~|~i\in I)$ for the topology $J$, the family
$(uX_i\to uX~|~i\in I)$ covers $uX$ for the topology $J'$.
\item[(C1)]
$u$ is continuous for the topologies $J$ and $J'$.
\item[(C2)]
Condition (C0) holds, and for every $X\in\Ob(\cC)$, every
covering family $(X_i\to X~|~i\in I)$ for the topology $J$
admits a refinement $(X'_i\to X~|~i\in I')$ that still covers
$X$, such that the fibre products $X'_i\times_XX'_j$ and
$X'_i\times_XX'_j\times_XX'_k$ are representable in $\cC$
for every $i,j,k\in I$ and $u$ commutes with these fibre products.
\item[(C3)]
Condition (C0) holds, and there exists a universe $\sV$ containing
$\sU$ such that $\cC$ and $\cC'$ are $\sV$-small and the functor
$u^a_{\sV!}:\cC^\wedge_\sV\to C'^\sim_\sV$ commutes with fibre products.
\item[(C4)]
$u:C'\to C$ is a morphism of sites.
\end{itemize}

\begin{remark}\label{rem_conditions-C0--C4}
(i)\ \
Recall that (C1) means that $u^\wedge$ transforms sheaves on
$C'$ into sheaves on $C$.  Likewise, according to lemma
\ref{lem_crit-continuity}, condition (C0) holds if and only
if $u^\wedge$ transforms separated presheaves on $C'$ into
separated presheaves on $C$.

(ii)\ \
We have (C4)$\Rightarrow$(C3)$\Rightarrow$(C1)$\Rightarrow$(C0),
by lemma \ref{lem_crit-continuity}. By inspecting the proof of
lemma \ref{lem_crit-continuity} we also easily see that
(C2)$\Rightarrow$(C1).

(iii)\ \
In the situation of \eqref{subsec_weak-strong-continuity},
suppose that $C$ is a lex-site and $u$ is left exact.
Then example \ref{ex_simple-case} says that $u$ fulfills
condition (C4) if and only if it fullfils condition (C0).
\end{remark}

\sset\subsubsection{}\label{subsec_pull-back-presh-of-cart}
In the situation of \eqref{subsec_weak-strong-continuity},
let $\cA'\to\cC'$ be a fibration, and $X\in\Ob(\cC)$. Set
$\cA:=\Fib(u)^*\cA'\to\cC$ and let $\pi:\cA\to\cA'$ be
the natural projection. For every cartesian section
$\sigma\in\cA'(uX)$ there exists a unique cartesian
section
$$
u^*_{|X}(\sigma)\in\cA(X)
\qquad\text{such that}\qquad
\pi\circ u^*_{|X}(\sigma)=\sigma\circ u_{|X}.
$$
Also, for every pair of cartesian sections
$\sigma,\tau\in\cA'(uX)$ and every natural $\cC$-transformation
$\beta:\sigma\Rightarrow\tau$ there exists a unique natural
$\cC$-transformation
$$
u^*_{|X}(\beta):u^*_{|X}(\sigma)\Rightarrow u^*_{|X}(\tau)
\qquad\text{such that}\qquad
\pi*u^*_{|X}(\beta)=\beta*u_{|X}.
$$
This characterization easily implies that the rules :
$\sigma\mapsto u^*_{|X}$ and $\beta\mapsto u^*_{|X}(\beta)$
for every such $\sigma$ and $\beta$, define a functor
$u^*_{|X}:\cA'(uX)\to\cA(X)$ fitting into the commutative diagram :
$$
\xymatrix{
\cA'(uX) \ar[r]^-{u^*_{|X}} \ar[d]_{\sev^{\cA'}_{uX}} &
\cA(X) \ar[d]^{\sev^\cA_X} \\
\cA'_{uX} \ar[r] & \cA_X
}$$
whose bottom horizontal arrow is the isomorphism of
categories induced by $\pi$, and whose vertical arrows
are the evaluation functors of \eqref{subsec_eval-functor}.
Since the latter are equivalences, it follows that the same
holds for $u^*_{|X}$. Notice that for every $(f:Y\to X)\in\Ob(\cC)$
and every cartesian section $\sigma\in\cA'(uX)$ we have :
$$
u^*_{|Y}(\sigma\circ(uf)_*)=u^*_{|X}(\sigma)\circ f_*
$$
(detail left to the reader). Hence, for every pair of
cartesian sections $\sigma,\tau\in\cA'(uX)$ and
every $(f:Y\to X)\in\Ob(\cC/X)$ we deduce a bijection
$$
\sCart(\sigma,\tau)(uf)\isom\sCart(u^*_{|X}\sigma,u^*_{|X}\tau)(f)
\qquad
\beta\mapsto u^*_{|Y}\beta.
$$
Finally, it is easily seen that this system of
bijections amounts to an isomorphism of presheaves
$$
u^\wedge_{|X}\cCart(\sigma,\tau)\isom
\cCart(u^*_{|X}\sigma,u^*_{|X}\tau).
$$

\begin{proposition}\label{prop_C0-C4}
In the situation of \eqref{subsec_weak-strong-continuity},
let $F:\cE\to\cC'$ be any fibration. We have :
\begin{enumerate}
\item
If $\cE$ is $0$-separated, and {\em(C0)} holds, then
$\Fib(u)^*(\cE)$ is $0$-separated.
\item
If $\cE$ is $1$-separated, and {\em(C1)} holds, then
$\Fib(u)^*(\cE)$ is $1$-separated.
\item
If {\em(C2)} holds, then $u$ is a weak morphism of sites $C'\to C$.
\end{enumerate}
\end{proposition}
\begin{proof} To ease notation, set $\cF:=\Fib(u)^*(\cE)$,
and let $\pi:\cF\to\cE$ be the natural projection.

(i): Let $X\in\Ob(\cC)$ be any object, and
$\sigma,\tau\in\cF(X)$ two cartesian sections;
according to lemma \ref{lem_crit-0-1_separation} it suffices
to show that the presheaf $\cCart(\sigma,\tau)$ is separated
on the site $C/X$. However, the discussion of
\eqref{subsec_pull-back-presh-of-cart} yields cartesian
sections $\sigma',\tau'\in\cE(uX)$ with isomorphisms
$u^*_{|X}\sigma'\isom\sigma$ and $u^*_{|X}\tau'\isom\tau$,
whence isomorphisms of presheaves
$$
\cCart(\sigma,\tau)\isom\cCart(u^*_{|X}\sigma',u^*_{|X}\tau')
\isom u^\wedge_{|X}\cCart(\sigma',\tau').
$$
Since $\cE$ is $0$-separated, lemma \ref{lem_crit-0-1_separation}
says that $\cCart(\sigma',\tau')$ is a separated presheaf on
the site $C'/uX$. On the other hand, the explicit description
of the covering sieves of the sites $C/X$ and $C'/uX$ furnished
by \eqref{sec_Localization-topoi} easily implies that the
functor $u_{|X}$ also fulfills condition (C0), relative to these
sites (details left to the reader). Then the assertion follows
from remark \ref{rem_conditions-C0--C4}(i).

(ii) is similar : arguing as in the foregoing, we are reduced
to checking that $u^\wedge_{|X}(\cCart(\sigma',\tau'))$ is a sheaf
on $C/X$. However, $\cCart(\sigma',\tau')$ is a sheaf on $C'/uX$,
by lemma \ref{lem_crit-0-1_separation}, and $u_{|X}$ is continuous
on the sites $C/X$ and $C'/uX$, by proposition
\ref{prop_localize-continuity}, whence the assertion.

(iii): Suppose that (C2) holds, and that $\cE$ is a stack on $C'$;
let $X\in\Ob(\cC)$ be any object, $\cS\subset\cC/X$ a sieve
covering $X$, and $\cS'\subset\cS$ a refinement covering $X$,
generated by a family $f_\bullet:=(f_i:X'_i\to X~|~i\in I)$ with
the properties of condition (C2). We get restriction functors :
$$
\cF(X)\xrightarrow{\rho}\sCart_\cC(\cS,\cF)
\xrightarrow{\rho'}\sCart_\cC(\cS',\cF)
$$
and we need to prove that $\rho$ is an equivalence. We check
first that $\rho$ is fully faithful. To this aim, notice that
$\cF$ is $0$-separated by virtue of (i), hence $\rho'$ is
faithful (lemma \ref{lem_pair-of-sieves}), and we are then
easily reduced to showing that $\rho'\circ\rho$ is fully
faithful. Now, choose a cleavage $\sc$ for $\cE$, so that
$\sc\circ u^o$ is a cleavage for $\cF$. Then $\sCart_\cC(\cS',\cF)$
is equivalent to the category of descent data
$\Desc(\cF,f_\bullet,\sc\circ u^o)$ as in \eqref{subsec_descnt-data}.
However, $\pi$ induces an isomorphism of categories
$\cF_X\isom\cE_{uX}$, and since $u$ commutes with the fibre
products $X_i\times_XX_j$ and $X_i\times_XX_j\times_XX_k$, a
simple inspection shows that the natural functor
$$
\cF_X\to\Desc(\cF,f_\bullet,\sc\circ u^o)
$$
is naturally identified with the corresponding functor
$$
\cE_{uX}\to\Desc(\cE,u(f_\bullet),\sc).
$$
The latter is an equivalence, since $\cE$ is a stack. Thus
$\rho'\circ\rho$ is even an equivalence, and therefore $\rho$
is fully faithful, as stated. Since $X$ and $\cS$ are arbitrary,
we have thus proved that $\cF$ is $1$-separated; but then $\rho'$
is fully faithful, again by lemma \ref{lem_pair-of-sieves}, and
it follows easily that $\rho$ is an equivalence if and only if
the same holds for $\rho'\circ\rho$. For the latter, the assertion
has just been shown, so finally $\cF$ is a stack, as required.
\end{proof}

\begin{lemma}\label{lem_kayak-days}
In the situation of \eqref{subsec_weak-strong-continuity},
let $\phi:\cE_1\to\cE_2$ be a $\cC'$-cartesian functor between
two fibrations over $\cC'$. Suppose that {\em(C0)} holds, and
moreover that for every $X\in\Ob(\cC')$ there exists a covering
family $(uY_j\to X~|~j\in I)$ for the topology $J'$. Let also
$i\in\{0,1,2\}$ such that $\Fib(u)^*(\phi)$ is $i$-covering for
the topology $J$. Then $\phi$ is $i$-covering for the topology $J'$.
\end{lemma}
\begin{proof} In view of lemma \ref{lem_reduce-to-split-fibs},
we may replace $\phi$ by $\sC(\phi)$, and reduce to the case
where $\cE_1$ and $\cE_2$ are split $\cC'$-fibrations, and
$\phi$ is a split cartesian functor. For $t=1,2$ let
$\sc_t:\cC'\to\bCat$ be the strict pseudo-functor associated
with the split cleavage of $\cE_t$. Suppose first that $i=0$, and
let $X\in\Ob(\cC')$ and $E\in\Ob(\cE_{2,X})$. By assumption, we
may find a covering family $(f_j:uY_j\to X~|~j\in I)$ and we set
$E_j:=\sc_{2,f_j}E\in\Ob(\cE_{2,uY_j})$ for every $j\in I$. Since
$\Fib(u)^*(\phi)$ is $0$-covering, we may find for every $j\in I$
a covering family $(f_{jj'}:Y_{jj'}\to Y_j~|~j'\in\Lambda_j)$ and for
every $j'\in\Lambda_j$ an object $E'_{jj'}\in\Ob(\cE_{1,uY_{jj'}})$
with an isomorphism $\phi E'_{jj'}\isom\sc_{2,u(f_{jj'})}E_j$. Set
$h_{jj'}:=f_j\circ u(f_{jj'}):uY_{jj'}\to X$ for every $j\in I$ and
$j'\in\Lambda_j$; since (C0) holds for $u$, the family
$(h_{jj'}~|~j\in I,\ j'\in\Lambda_j)$ covers $X$, and
$\phi E'_{jj'}\isom\sc_{2,h_{jj'}}E$ for every $j\in I$ and
$j'\in\Lambda_j$. This shows that $\phi$ is $0$-covering.

Next, suppose that $i=1$, and let $X\in\Ob(\cC')$ and
$g:\phi E_1\to\phi E_2$ a morphism in $\cE_{2,X}$. We pick
a covering family $(f_j:uY_j\to X~|~j\in I)$, and set
$E_{tj}:=\sc_{1,f_j}E_t$ for $t=1,2$ and every $j\in I$;
also, let $g_j:=\sc_{2,f_j}(g):\phi E_{1j}\to\phi E_{2j}$ for
every $j\in J$. Since $\Fib(u)^*(\phi)$ is $1$-covering,
for every $j\in I$ there exist a covering family
$(f_{jj'}:Y_{jj'}\to Y_j~|~j'\in\Lambda_j)$ and for every
$j'\in\Lambda_j$ a morphism
$g_{jj'}:\sc_{1,u(f_{jj'})}E_{1j}\to\sc_{1,u(f_{jj'})}E_{2j}$ with
$\phi(g_{jj'})=\sc_{2,u(f_{jj'})}(g_j)=\sc_{2,u(f_{jj'})\circ f_j}(g)$.
Arguing as in the foregoing, we easily deduce that $\phi$
is $1$-covering (details left to the reader).

Lastly, suppose that $i=2$, and let $X\in\Ob(\cC')$ and
$g_1,g_2:E_1\to E_2$ two morphisms in $\cE_{1,X}$ such that
$\phi(g_1)=\phi(g_2)$. Pick again a covering family
$(f_j:uY_j\to X~|~j\in I)$ and set $g_{tj}:=\sc_{1,f_j}(g_t)$
for $t=1,2$ and every $j\in I$. Then $\phi(g_{1j})=\phi(g_{2j})$
for every $j\in I$, and since $\Fib(u)^*(\phi)$ is $2$-covering,
there exists for every such $j$ a covering family
$(f_{jj'}:Y_{jj'}\to Y_j~|~j'\in\Lambda_j)$ such that
$\sc_{1,u(f_{jj'})}(g_{1j})=\sc_{1,u(f_{jj'})}(g_{2j})$ for every
$j'\in\Lambda_j$. Arguing as in the foregoing, it follows
easily that $\phi$ is $2$-covering.
\end{proof}

\sset\subsubsection{}\label{subsec_extend-stacks}
Let $C:=(\cC,J)$ be a small site and define the site
$C^\wedge_\sU:=(\cC^\wedge_\sU,J^\wedge)$ as in remark
\ref{rem_topol-on-presheaves}(ii). By theorem
\ref{th_canon-topos}(i), the Yoneda embedding is a
cocontinuous morphism of sites $h_\cC:C^\wedge_\sU\to C$.
Thus, for every universe $\sV$ containing $\sU$ and
every $i$-separated $\sV$-prestack $\phi:\cE\to\cC$
on $\cC$, the fibration
$$
\sV\tdu\Fib(h_\cC)_*(\cE)\to\cC^\wedge_\sU
$$
is an $i$-separated prestack for the topology $J^\wedge$
(proposition \ref{prop_cocont-and-stacks}(ii)). Moreover,
since $h_\cC$ is fully faithful, remark
\ref{rem_explicit-Fib_*}(iv) and corollary
\ref{cor_fully-faith-2-adjoint} imply that the counit of
$2$-adjunction
$$
\sV\tdu\Fib(h_\cC)^*\circ\sV\tdu\Fib(h_\cC)_*(\cE)\to\cE
$$
is an equivalence of categories. Thus,
$\sV\tdu\Fib(h_\cC)_*(\cE)$ is a natural extension of $\cE$
to the site $C^\wedge$ which contains $C$ as a full subcategory
with the topology induced from $J^\wedge$.

\begin{remark}\label{rem_equiv-of-stack-cats}
(i)\ \
In the situation of \eqref{subsec_extend-stacks}, we can
describe more explicitly the fibration $\sV\tdu\Fib(h_\cC)_*(\cE)$,
as follows. We consider the functor
$$
\sc^\triangledown_\phi:(\cC^\wedge_\sU)^o\to\sV\tdu\bCat
\qquad
F\mapsto\sCart_\cC(\cFib(F),\cE)
$$
which assigns to every morphism of presheaves $\beta:F\to F'$,
the functor
$\sCart_\cC(\cFib(\beta),\cE):\sc^\triangledown_\phi(F')\to
\sc^\triangledown_\phi(F)$. Then, combining remark
\ref{rem_explicit-Fib_*}(ii) and example
\ref{ex_fibred-cats-II}(iii) we get a natural equivalence
of $\cC^\wedge_\sU$-categories :
$$
\sV\tdu\Fib(h_\cC)_*(\cE)\isom\cFib(\sc^\triangledown_\phi).
$$

(ii)\ \
Let $i:C^\sim_\sU\to\cC^\wedge_\sU$ be the inclusion functor, and
$(-)^a:\cC^\wedge_\sU\to C^\sim_\sU$ its left adjoint. We have
a pseudo-natural equivalence of pseudo-functors :
$$
\sV\tdu\Fib((-)^a)_*\isom\sV\tdu\Fib(i)^*.
$$
Indeed, let $\cE$ be any fibration on $\cC^\wedge_\sU$; in view
of remark \ref{rem_explicit-Fib_*}(ii), the fibration
$\sV\tdu\Fib((-)^a)_*(\cE)$ is naturally equivalent to the
fibration associated with the strict pseudo-functor
$$
F\mapsto\sCart_{\cC^\wedge_\sU}((-)^a\cC^\wedge_\sU/F,\cE)
\qquad
\text{for every $F\in\Ob(C^\sim_\sU)$}.
$$
On the other hand, $\sC(\Fib(i)^*(\cE))$ is the fibration associated
with the strict pseudo-functor
$$
F\mapsto\sCart_{\cC^\wedge_\sU}(\cC^\wedge_\sU/iF,\cE)
\qquad
\text{for every $F\in\Ob(C^\sim_\sU)$}.
$$
It suffices then to notice that $(-)^a\cC^\wedge_\sU/F=\cC^\wedge_\sU/iF$
for every sheaf $F$ on the site $C$.
\end{remark}

\begin{theorem}\label{th_2-equiv-for-cats-of-stacks}
Let $C:=(\cC,J)$ be a small site, $h_\cC:\cC\to\cC^\wedge_\sU$
the Yoneda embedding and $i:C^\sim_\sU\to\cC^\wedge_\sU$ the
inclusion functor. For every universe $\sV$ containing
$\sU$ we have :
\begin{enumerate}
\item
The strict pseudo-functor $\sV\tdu\Fib(h_\cC)^*$ restricts
to a $2$-equivalence ;
$$
\sV\tdu\Stack(C^\wedge_\sU)\isom\sV\tdu\Stack(C).
$$
\item
Likewise, $i$  and its left adjoint
$(-)^a:\cC^\wedge_\sU\to C^\sim_\sU$ induce $2$-equivalences :
$$
\xymatrix{
\sV\tdu\Stack(C^\wedge_\sU) \ar@<.5ex>[rr]^-{\sV\tdu\Fib(i)^*} & &
\sV\tdu\Stack(\Can(C^\sim_\sU)).
\ar@<.5ex>[ll]^-{\sV\tdu\Fib((-)^a)^*}
}$$
\end{enumerate}
\end{theorem}
\begin{proof}(i): Let $\eta$ and $\eps$ be the unit and counit
of a $2$-adjunction for the pair of pseudo-functors
$(\sV\tdu\Fib(h_\cC)^*,\sV\tdu\Fib(h_\cC)_*)$. Since $h_\cC$
is a cocontinuous morphism of sites $C^\wedge\to C$ (theorem
\ref{th_canon-topos}(i)), we know already that
$\sV\tdu\Fib(h_\cC)_*$ sends $\sV$-stacks on $C$ to $\sV$-stacks
on $C^\wedge$ (proposition \ref{prop_cocont-and-stacks}(ii));
according to corollary \ref{cor_characterize-2-eq}(i), it then
suffices to show that for every $\sV$-stack $\cE$ on $C^\wedge$
and $\cE'$ on $C$, the fibration $\sV\tdu\Fib(h_\cC)^*(\cE)$ is a
$\sV$-stack on $C$, and $\eta_\cE$ and $\eps_{\cE'}$ are equivalences.
The assertion for $\eps_{\cE'}$ has already been noticed in
\eqref{subsec_extend-stacks}. Next, to ease notation, let us
drop the prefix $\sV$, and set $\cF:=\Fib(h_\cC)^*(\cE)$. Let
also $i_\cF:\cF\to\cF^a$ be the natural morphism of fibrations
on $\cC$, and denote by
$$
\beta_\cE:\cE\to\Fib(h_\cC)_*(\cF^a)
$$
the composition of $\eta_\cE$ and $\Fib(h_\cC)_*(i_\cF)$.
Then $\eps_{\cF^a}\circ(\Fib(h_\cC)^*\beta_\cE)$ is isomorphic to
$$
\eps_{\cF^a}\circ(\Fib(h_\cC)^*\circ\Fib(h_\cC)_*(i_\cF))\circ
(\Fib(h_\cC)^*\eta_\cE)
$$
which is in turn isomorphic to
$i_\cF\circ\eps_\cF\circ(\Fib(h_\cC)^*\eta_\cE)$, and the latter is
isomorphic to $i_\cF$, by virtue of the triangular modifications
associated with the pair $(\eta,\eps)$ (theorem
\ref{th_2-adjunction}(i)). Since $\eps_{\cF^a}$ is an equivalence,
and $i_\cF$ is $j$-covering for $j=0,1,2$, it follows that
$\Fib(h_\cC)^*(\beta_\cE)$ is $j$-covering for $j=0,1,2$. Since
the functor $h_\cC$ verifies condition (C0) of
\eqref{subsec_weak-strong-continuity}, lemma \ref{lem_kayak-days}
implies that $\beta_\cE$ is also $j$-covering for $j=0,1,2$. But
$\Fib(h_\cC)_*(\cE)$ is a stack (proposition
\ref{prop_cocont-and-stacks}(ii)), hence $\beta_\cE$ is an
equivalence (proposition \ref{prop_faith-separ-cover}).
Therefore, $\Fib(h_\cC)^*(\beta)$ is an equivalence as well,
and by the foregoing, the same then holds for $i_\cF$. The
latter means that $\Fib(h_\cC)^*(\cE)$ is a stack, and it also
follows that $\eta_\cE$ is an equivalence, as sought.

(ii): By theorem \ref{th_canon-topos}(iii), the functor $(-)^a$
is cocontinuous for the topologies $J^\wedge$ and $\Can_{C^\sim}$;
taking into account remark \ref{rem_equiv-of-stack-cats}(ii)
and proposition \ref{prop_cocont-and-stacks}(ii), it follows
that $i$ is a weak morphism of sites $C^\wedge_\sU\to\Can(C^\sim_\sU)$.
Also, $(-)^a$ is a weak morphism of sites
$\Can(C^\sim_\sU)\to C^\wedge_\sU$, due to proposition
\ref{prop_C0-C4}(iii).
Moreover, since $(-)^a\circ i$ is isomorphic to $\one_{C^\sim}$,
we have a pseudo-natural isomorphism of pseudo-functors:
$$
\sV\tdu\Fib(i)^*\circ\sV\tdu\Fib((-)^a)^*\isom
\one_{\sV\tdu\Stack(\Can(C^\sim))}.
$$
Next, let $\cE$ be any stack on $C^\wedge$ with $\sV$-small
fibres; recall that $\cE$ is naturally equivalent to the split
fibration associated with the strict pseudo-functor
$\sCart_{\cC^\wedge}(\cC^\wedge/-,\cE)$. Likewise,
$\sV\tdu\Fib((-)^a)^*\circ\sV\tdu\Fib(i)^*(\cE)$ is naturally
equivalent to the split fibration associated with the strict
pseudo-functor
$$
\cC^\wedge\to\sV\tdu\bCat
\qquad
F\mapsto\sCart_{\cC^\wedge}(\cC^\wedge/F^a,\cE).
$$
By virtue of corollary \ref{cor_characterize-2-eq}(i), we are
thus reduced to checking that for every $F\in\Ob(\cC^\wedge)$,
the natural morphism $j_F:F\to F^a$ induces an equivalence of
categories $\sCart_{\cC^\wedge}(\cC^\wedge/F^a,\cE)\isom
\sCart_{\cC^\wedge}(\cC^\wedge/F,\cE)$. To this aim, we notice :

\begin{claim}\label{cl_bicover-gives-eq}
For every stack $\cE$ on $C^\wedge$ and every bicovering morphism
$f:G\to G'$ of $(\cC^\wedge)^\wedge_\sV$ (for the topology $J^\wedge$
of $\cC^\wedge$), the functor $\sCart_{\cC^\wedge}(\cFib(f),\cE)$ is
an equivalence.
\end{claim}
\begin{pfclaim} Indeed, according to theorem \ref{th_stackfication}
we have an essentially commutative diagram :
$$
\xymatrix{ \sCart_\cC(\cFib(G'^a),\cE)
\ar[rrr]^-{\sCart_\cC(\cFib(f^a),\cE)} \ar[d] & & &
\sCart_\cC(\cFib(G^a),\cE) \ar[d] \\
\sCart_\cC(\cFib(G'),\cE) \ar[rrr]^-{\sCart_\cC(\cFib(f),\cE)}
& & & \sCart_\cC(\cFib(G),\cE)
}$$
whose vertical arrows are equivalences. By assumption,
$f^a:G^a\to G'^a$ is an isomorphism of sheaves on
$(C^\wedge)^\sim_\sV$, hence the top horizontal arrow is
an isomorphism of categories, and therefore the bottom
horizontal arrow is an equivalence, as claimed.
\end{pfclaim}

In view of claim \ref{cl_bicover-gives-eq}, and recalling
the natural identifications
$$
\cC^\wedge/F\isom\cFib(h_F)
\qquad
\cC^\wedge/F^a\isom\cFib(h_{F^a})
$$
(example \ref{ex_fibred-cats-II}(i)), we are then further
reduced to showing :

\begin{claim}
$j_F$ induces a bicovering morphism $h_{j_F}:h_F\to h_{F^a}$
of $(\cC^\wedge)^\wedge_\sV$.
\end{claim}
\begin{pfclaim}[] Consider any $\sV$-presheaf $G$ on
$\cC^\wedge$, denote by $G^a\in(C^\wedge)^\sim_\sV$ the sheaf
associated with $G$, and let $j_G:G\to G^a$ be the
natural morphism; since $h_\cC$ is continuous and cocontinuous
for the topologies $J$ and $J^\wedge$ (theorem
\ref{th_canon-topos}(i)), lemma \ref{lem_improve}(ii) implies
that $h^\wedge_\cC(G^a)$ is a $\sV$-sheaf on $C$, and the morphism
of $\sV$-presheaves $h^\wedge_\cC(G)\to h^\wedge_\cC(G^a)$ on $\cC$
is bicovering. Taking $G:=h_F$, we get $h^\wedge_\cC(h_F)=F$, whence
an isomorphism $t:h^\wedge_\cC(h^a_F)\isom F^a$ that identifies
$h^\wedge_\cC(j_G)$ with $j_F$. Likewise, we have
$h^\wedge_\cC(h_{F^a})=F^a$, and since $h^\wedge_\cC$ restricts to
an equivalence $(C^\wedge)^\sim_\sV\isom C^\sim_\sV$ (theorem
\ref{th_canon-topos}(iii)), there follows a unique isomorphism
$s:h^a_F\isom h_{F^a}$ such that $h_\cC^\wedge(s)=t$. Summing up,
we deduce that $h^\wedge_\cC(s\circ j_G)=j_F=h^\wedge_\cC(h_{j_F})$;
since the functor $h^\wedge_{\cC^\wedge}:\cC^\wedge\to(\cC^\wedge)^\wedge_\sV$
is fully faithful, we get $s\circ j_G=h_{j_F}$, and since $s$ is an
isomorphism, the claim follows.
\end{pfclaim}
\end{proof}

\begin{corollary}\label{cor_another-condition}
In the situation of \eqref{subsec_weak-strong-continuity},
suppose that condition {\em (C3)} holds for some universe $\sV$.
Then $u$ is a weak morphism of sites $C'\to C$.
\end{corollary}
\begin{proof} By remark \ref{rem_change-universe}(iii) we may
assume that also $\cC^\wedge$ and $\cC'^\wedge$ are $\sV$-small.
We consider the essentially commutative diagram :
$$
\xymatrix{ \cC \ar[rr]^-u \ar[d]_{h_\cC} & &
\cC' \ar[d]^{h^a_{C'}} \\
\cC^\wedge_\sV \ar[rr]^-{u^a_{\sV!}} & & C'^\sim_\sV
}$$
and we endow $\cC^\wedge_\sV$ with the topology $J^\wedge$ as
in remark \ref{rem_topol-on-presheaves}, and $C'^\sim_\sV$
with its canonical topology. We notice that, with these
topologies, the functor $u^a_{\sV!}$ satisfies condition
(C2) : indeed, condition (C0) holds for $u^a_{\sV!}$ by
definition of the topology $J^\wedge$; also, all fibre
products of $\cC^\wedge_\sV$ are representable, and by
assumption $u^a_{\sV!}$ commutes with fibre products. Let
$\sV'$ be a universe containing $\sV$ and such that
$\sC^\wedge_\sV$ and $C'^\sim_\sV$ are $\sV'$-small; by
the foregoing case, we deduce that $\sV'\tdu\Fib(u^a_{\sV!})^*$
sends stacks on $\Can(C'^\sim_\sV)$ to stacks on
$(\cC^\wedge_\sV,J^\wedge)$. On the other hand, theorem
\ref{th_2-equiv-for-cats-of-stacks} implies that
$\sV'\tdu\Fib(h^a_{C'})^*$ restricts to a $2$-equivalence from
the $2$-category of $\sV'$-stacks on $\Can(C'^\sim_\sV)$ to
that of $\sV'$-stacks on $(\cC',J')$, so we may assume that
$\cE=\sV'\tdu\Fib(h^a_{C'})^*(\cE')$ for some $\sV'$-stack
on $\Can(C'^\sim_\sV)$, and we set
$\cE'':=\sV'\tdu\Fib(u^a_{\sV!})^*(\cE')$. Then $\Fib(u)^*(\cE)$
is equivalent to $\sV'\tdu\Fib(h_\cC)^*(\cE'')$, and the latter
is a stack, by theorem \ref{th_2-equiv-for-cats-of-stacks}(i).
\end{proof}

\begin{example}\label{ex_localize-is-weak}
Let $\sV$ be a universe, $C:=(\cC,J)$ a $\sV$-small
site, and $X\in\Ob(\cC)$. Define the site $C/X$ as in
\eqref{sec_Localization-topoi}, and recall that the source
functor $\ss_X:\cC/X\to\cC$ is continuous and cocontinuous
for the topologies of $C$ and $C/X$. Moreover, the functor
$\ss_{\sV!}$ commutes with fibre products, by proposition
\ref{prop_in-the-same-vein}(vi.c), hence condition (C3) holds
for $\ss_X$, and by corollary \ref{cor_another-condition} we
conclude that $\ss_X$ is a weak morphism of sites $C\to C/X$.
\end{example}

\begin{proposition}\label{prop_equiv-and-2-equiv}
Let $C:=(\cC,J)$ and $C':=(\cC',J')$ be two $\sU$-sites, and
$u:C\to C'$ a morphism of sites such that :
\begin{enumerate}
\alphaenu
\item
$u$ induces an equivalence of categories
$\tilde u_*:C^\sim\isom C'^\sim$.
\item
The subset $u(\Ob(\cC'))\subset\Ob(\cC)$ is a topologically
generating family for the site $C$.
\end{enumerate}
Then, for every universe $\sV$ such that $C$ is a $\sV$-site,
the following holds :
\begin{enumerate}
\item
$u$ induces an equivalence of categories
$\tilde u_{\sV*}:C^\sim_\sV\isom C'^\sim_\sV$.
\item
$u$ induces a $2$-equivalence of\/ $2$-categories
$\sV\tdu\Fib(u)^*:\sV\tdu\Stack(C)\isom\sV\tdu\Stack(C')$.
\end{enumerate}
\end{proposition}
\begin{proof}(i): To ease notation, set $v:=\tilde u_*$.
Recall that $\Can(C^\sim)$ and $\Can(C'^\sim)$ are isomorphic
to $\sU$-sites (remark \ref{rem_summarized}(iv)). By corollary
\ref{cor_two-U-sites}(i.b) and theorem \ref{th_canon-topos}(iii)
there follows, for every universe $\sV$ containing $\sU$, an
essentially commutative diagram
$$
\xymatrix{
\Can(C^\sim)^\sim_\sV \ar[rr]^-{\tilde v_{\sV*}} \ar[d] & &
\Can(C'^\sim)^\sim_\sV \ar[d] \\
C^\sim_\sV \ar[rr]^-{\tilde u_{\sV*}} & & C'^\sim_\sV 
}$$
whose vertical arrows are equivalences. Since $v$ is an
equivalence, the same holds for $\tilde v_*$, and then
also for $\tilde u_*$. Next, let $\sV$ be any universe
such that $C$ is a $\sV$-site, and pick another universe
$\sV'$ containing $\sV$ and $\sU$. We deduce a commutative
diagram
$$
\xymatrix{ C^\sim_\sV \ar[rr]^-{\tilde u_{\sV*}} \ar[d] & &
C'^\sim_\sV \ar[d] \\
C^\sim_{\sV'} \ar[rr]^-{\tilde u_{\sV'*}} & & C'^\sim_{\sV'} 
}$$
whose vertical arrows are the fully faithful inclusion functors
and whose bottom horizontal arrow is an equivalence, by the
foregoing case. It follows already that $\tilde u_{\sV*}$ is
fully faithful. We may now argue as in the proof of theorem
\ref{th_gener-Beilinson} : let $F$ be any $\sV$-sheaf on $C'$;
we know already that there exists a $\sV'$-sheaf $G$ on $C$
such that $\tilde u_{*\sV'}G$ is isomorphic to $F$, and we need
to check that $GX$ is essentially $\sV$-small for every
$X\in\Ob(\cC)$. But by condition (b) we may find a covering
family $(uY_i\to X~|~i\in I)$, and by claim \ref{cl_reduce-size}
we may assume that $I$ is small. Then the induced map
$GX\to\prod_{i\in I}G(uY_i)$ is injective, and by assumption
$G(uY_i)=FY_i$ is $\sV$-small for every $i\in I$, whence the
contention.

(ii): Suppose first that $\sV$ contains $\sU$ and both
$\cC^\wedge$ and $\cC'^\wedge$ are $\sV$-small; notice that
$u$ fulfills condition (C3) of
\eqref{subsec_weak-strong-continuity} (remark
\ref{rem_conditions-C0--C4}(ii)), and therefore
$\sV\tdu\Fib(u)^*$ sends $\sV$-stacks on $C$ to $\sV$-stacks
on $C'$ (corollary \ref{cor_another-condition}). Likewise,
$\sV\tdu\Fib(v)^*$ sends $\sV$-stacks on $\Can(C^\sim)$ to
$\sV$-stacks on $\Can(C'^\sim)$. By corollary
\ref{cor_two-U-sites}(i.b) and theorem
\ref{th_2-equiv-for-cats-of-stacks}(ii) we then get a
pseudo-commutative diagram
$$
\xymatrix{
\sV\tdu\Stack(\Can(C^\sim)) \ar[rr]^-{\sV\tdu\Fib(v)^*} \ar[d] & &
\sV\tdu\Stack(\Can(C'^\sim)) \ar[d] \\
\sV\tdu\Stack(C) \ar[rr]^-{\sV\tdu\Fib(u)^*} & & \sV\tdu\Stack(C')
}$$
whose vertical arrows are $2$-equivalences. Since $v$ is an
equivalence, $\sV\tdu\Fib(v)^*$ is also a $2$-equivalence, and
then the same holds for $\sV\tdu\Fib(u)^*$. Lastly, let $\sV$
be any universe such that $C$ is a $\sV$-small site, and pick
a universe $\sV'$ containing $\sV$ and $\sU$, and such that
$\cC^\wedge$ and $\sC'^\wedge$ are $\sV'$-small; we deduce a
commutative diagram
$$
\xymatrix{ \sV\tdu\Stack(C) \ar[rr]^-{\sV\tdu\Fib(u)^*} \ar[d] & &
\sV\tdu\Stack(C') \ar[d] \\
\sV'\tdu\Stack(C) \ar[rr]^-{\sV'\tdu\Fib(u)^*} & & \sV'\tdu\Stack(C') 
}$$
whose vertical arrows are the fully faithful inclusion strict
pseudo-functors, and whose bottom horizontal arrow is a
$2$-equivalence, by the foregoing case. It follows already
that $\sV\tdu\Fib(u)^*$ is fully faithful. We remark :

\begin{claim}\label{cl_small-stacks}
Let $T\subset\Ob(\cC)$ be a $\sV$-small topologically generating
family of the site $C$, and $F:\cE\to\cC$ a $\sV'$-stack on $C$
such that $F^{-1}X$ is a $\sV$-small category for every $X\in T$.
Then $\cE$ is equivalent to a $\sV$-stack on $C$.
\end{claim}
\begin{pfclaim} Denote by $\cT\subset\cC$ the full subcategory
with $\Ob(\cT)=T$, let $j:\cT\to\cC$ be the inclusion functor,
and endow $\cT$ with the topology $J_\cT$ induced by $J$ via
the functor $j$. Let also
$\eta_\cE:\cE\to\cF:=\sV'\tdu\Fib(j)_*\circ\sV'\tdu\Fib(j)^*(\cE)$
be the unit of adjunction. By proposition
\ref{prop_comparis-lemma}(ii), the functor $j$ is a morphism of
$\sV'$-sites $C\to(\cT,J_\cT)$, hence $\cE':=\sV'\tdu\Fib(j)^*(\cE)$
is a $\sV'$-stack on $(\cT,J_\cT)$ (remark
\ref{rem_conditions-C0--C4}(ii) and corollary
\ref{cor_another-condition}), and by assumption its fibres are
$\sV$-small, {\em i.e.} it is a $\sV$-stack on $(\cT,J_\cT)$.
Moreover, since $j$ is cocontinuous for the topologies $J$ and
$J_\cT$ (proposition \ref{prop_comparis-lemma}(i)), we also see
that $\cF$ is a $\sV'$-stack on $C$ (proposition
\ref{prop_cocont-and-stacks}(ii)). Furthermore, since $j$ is
fully faithful, the same holds for $\sV'\tdu\Fib(j)_*$ (remark
\ref{rem_explicit-Fib_*}(iv)), so the counit of adjunction
$\eps_{\cE'}:\sV'\tdu\Fib(j)^*(\cF)\to\cE'$ is an equivalence.
By virtue of the triangular modifications for the pair
$(\eta,\eps)$ (theorem \ref{th_2-adjunction}(i)), it follows
that $\sV'\tdu\Fib(j)^*(\eta_\cE)$ is an equivalence as well,
and then it is $i$-covering for $i=0,1,2$ (proposition
\ref{prop_faith-separ-cover}). By invoking lemma \ref{lem_kayak-days}
we conclude finally that $\eta_\cE$ is $i$-covering for $i=0,1,2$,
so it is an equivalence, again by proposition
\ref{prop_faith-separ-cover}. Summing up, we are reduced to
checking that if $\cE'$ is a $\sV$-stack on $(\cT,J_\cT)$, then
$\cF':=\sV'\tdu\Fib(j)_*(\cE')$ is a $\sV$-stack on $C$. To this
aim, pick any cleavage $\sc:\cT^o\to\sV\tdu\bCat$ for $\cE'$; by
construction, the fibre of $\cF'$ over any $X\in\Ob(\cC)$ is the
$2$-limit of the pseudo-functor
$\sc\circ\st_{X^o}:X^o/j^o\cT^o\to\sV\tdu\bCat$, and it suffices to
remark that $X^o/j^o\cT^o$ is a $\sV$-small category.
\end{pfclaim}

Now, by assumption $\cC$ admits a $\sV$-small topologically
generating family $T$, and by assumption (b), for every $X\in T$
we may find a covering family $(uY_i\to X~|~i\in I_X)$. By
claim \ref{cl_reduce-size}, for every $X\in T$ there exists
then a $\sV$-small subset $I'_X\subset I_X$ such that
$(uY_i\to X~|~i\in I'_X)$ is still a covering family; set
$T':=\{uY_i~|~X\in T,\ i\in I_X\}$. Clearly $T'$ is a $\sV$-small
topologically generating family for $C$. Lastly, let $\cE'$
be any $\sV$-stack on $C'$; we know already that there exists
a $\sV'$-stack $\cE$ on $C$ such that $\sV'\tdu\Fib(u)^*(\cE)$
is equivalent to $\cE'$. But the latter condition implies that
the fibre of $\cE$ over each $X\in T'$ is a $\sV$-small
category; by claim \ref{cl_small-stacks} we deduce that
$\cE$ is equivalent to a $\sV$-stack on $C$, and this shows
that $\sV\tdu\Fib(u)^*$ is a $2$-equivalence, as stated.
\end{proof}

\begin{corollary}\label{cor_stackificate-U-sites}
For every $\sU$-site $C:=(\cC,J)$ the inclusion pseudo-functor
$\Stack(C)\to\Fib(\cC)$ admits a left $2$-adjoint
$$
\Fib(\cC)\to\Stack(C)
\qquad
\cE\mapsto\cE^a.
$$
\end{corollary}
\begin{proof} Let $T\subset\Ob(\cC)$ be a small topologically
generating family, $\cT\subset\cC$ the full subcategory
with $\Ob(\cT)=T$, and $j:\cT\to\cC$ the inclusion functor.
Endow $\cT$ with the topology $J_\cT$ induced by $J$ via $j$.
Let also $\cE\to\cC$ be any fibration with small fibres, and
set $\cE':=\Fib(j)^*(\cE)$. Let $\cE'^a$ be the stack associated
with $\cE'$ (theorem \ref{th_stackfication}), and denote by
$i:\cE'\to\cE'^a$ the natural morphism of fibrations over $\cT$.
Choose any universe $\sV$ containing $\sU$ and
such that $\cC$ is $\sV$-small, and consider the composition
$$
\beta_\cE:\cE\xrightarrow{\ \eta_\cE\ }
\sV\tdu\Fib(j)_*(\cE')\xrightarrow{\ \sV\tdu\Fib(j)_*(i)\ }
\cE^*:=\sV\tdu\Fib(j)_*(\cE'^a)
$$
where $\eta_\cE$ is the unit of adjunction. By arguing as in
the proof of claim \ref{cl_small-stacks}, we see that $\cE^*$
is a stack with small fibres on $C$ and that
$\sV\tdu\Fib(j)^*(\eta_\cE)$ is an equivalence of categories.
Moreover, since $j$ is fully faithful, the counit of adjunction
$\sV\tdu\Fib(j)^*\circ\sV\tdu\Fib(j)_*\Rightarrow
\one_{\sV\tdu\Fib(\cT)}$ is a pseudo-natural equivalence, so
we have an essentially commutative diagram of functors
$$
\xymatrix{ \sV\tdu\Fib(j)^*\circ\sV\tdu\Fib(j)_*(\cE')
\ar[rrrr]^-{\sV\tdu\Fib(j)^*\circ\sV\tdu\Fib(j)_*(i)}
\ar[d] & & & & \sV\tdu\Fib(j)^*(\cE^*) \ar[d] \\
\cE' \ar[rrrr]^-i & & & & \cE'^a
}$$
whose vertical arrows are equivalences of categories. Since
$i$ is $l$-covering for $l=0,1,2$ (proposition
\ref{prop_unit-sep-cov-faith}(i)), we conclude that the same
holds for $\sV\tdu\Fib(j)^*(\beta_\cE)$, and then also for
$\beta_\cE$, by virtue of lemma \ref{lem_kayak-days}. Let us
denote by $\cE^a_\sV$ the $\sV$-stack associated with $\cE$
(a priori, its fibres are only $\sV$-small), and
$i_\sV:\cE\to\cE^a_\sV$ the unit of adjunction; there exists
a morphism of fibrations $\omega_\cE:\cE^a_\sV\to\cE^*$,
pseudo-natural in $\cE$, such that the functor
$\omega_\cE\circ i_\sV$ is isomorphic to $\beta_\cE$. Since
$i_\sV$ is $l$-covering for $l=0,1,2$ (proposition
\ref{prop_unit-sep-cov-faith}(i)), it follows that the same
holds for $\omega_\cE$ (lemma \ref{lem_yoga-i-coverings}(ii)),
and then the latter is an equivalence of categories
(proposition \ref{prop_faith-separ-cover}). Summing up,
the rule : $\cE\mapsto\omega_\cE$ yields a pseudo-natural
equivalence of pseudo-functors, so the pseudo-functor
given by the rule : $\cE\mapsto\cE^*$ is the sought left
$2$-adjoint.
\end{proof}

\sset\subsubsection{}\label{subsec_St(u)_*}
Let now  $C:=(\cC,J)$ and $C':=(\cC',J')$ be two
$\sU$-sites, $u:C\to C'$ a weak morphism of sites
(definition \ref{def_weak-morph-of-sites}); then,
for every universe $\sV$ with $\sU\subset\sV$, the
pseudo-functor $\sV\tdu\St(u)_*$ admits a left $2$-adjoint
$$
\sV\tdu\St(u)^*:\sV\tdu\Stack(C')\to\sV\tdu\Stack(C)
$$
defined as follows. We choose a $\sU$-small topologically
generating subset $G\subset\Ob(\cC')$ and let $\cG\subset\cC'$
be the full subcategory with $\Ob(\cG)=G$. We endow $\cG$ with
the topology $J_\cG$ induced by $J'$ via the inclusion
functor $j:\cG\to\cC'$. Recall that $j$ is a morphism of
sites $C'\to(\cG,J_\cG)$ and
$\tilde\jmath_*:C'^\sim\isom(\cG,J_\cG)^\sim$ is an
equivalence (proposition \ref{prop_comparis-lemma}(ii)).
Then $\sV\tdu\St(j)_*$ is well defined (remark
\ref{rem_conditions-C0--C4}(ii)) and it is a $2$-equivalence
(proposition \ref{prop_equiv-and-2-equiv}(ii)), and
$\sV\tdu\St(u)^*$ is given by the composition :
$$
\sV\tdu\Stack(C')\xrightarrow{\sV\tdu\St(j)_*}\sV\tdu\Stack(\cG,J_\cG)
\xrightarrow{\sV\tdu\Fib(u\circ j)_!}\sV\tdu\Fib(\cC)\xrightarrow{(-)^a}
\sV\tdu\Stack(C)
$$
where $(-)^a$ is the pseudo-functor of corollary
\ref{cor_stackificate-U-sites}. It is then easily seen
that the resulting pseudo-functor is a left $2$-adjoint
for $\sV\tdu\St(u)_*$; especially, it is independent, up to
pseudo-natural equivalence, of the choice of $G$ (remark
\ref{rem_opposing-thumb}(iii)); the details shall be
left to the reader. As usual, we shall often omit mentioning
$\sV$, in case $\sV=\sU$.

\begin{example}\label{ex_when-St-upper-star-exists}
By remark \ref{rem_conditions-C0--C4}(ii), every
morphism of sites $u:C:=(\cC,J')\to C':=(\cC',J)$ fulfills
(C3), hence it is a weak morphism of sites (corollary
\ref{cor_another-condition}), and so it induces a pseudo-functor
$\sV\tdu\St(u)_*$. If $C$ and $C'$ are $\sU$-sites, we have
also the pseudo-functor $\sV\tdu\St(u)^*$.
\end{example}

\begin{proposition}\label{prop_Fib-lower-shriek-and-u}
Let $C:=(\cC,J)$ and $C':=(\cC',J')$ be two small sites,
$u:C\to C'$ a morphism of sites, $i\in\{0,1,2\}$, and
$f:\cE\to\cF$ an $i$-covering cartesian functor of\/
$\cC'$-fibrations. We have :
\begin{enumerate}
\item
The functor $\Fib(u)_!(f):\Fib(u)_!(\cE)\to\Fib(u)_!(\cF)$
is $i$-covering.
\item
If\/ $\cE$ and $\cF$ are stacks on $C'$, then
$\St(u)^*(f):\St(u)^*(\cE)\to\St(u)^*(\cF)$ is $i$-covering.
\end{enumerate}
\end{proposition}
\begin{proof} We show first the following special case :

\begin{claim}\label{cl_ok-for-lex-sites}
The proposition holds if $u$ is a morphism of lex-sites.
\end{claim}
\begin{pfclaim} By proposition \ref{prop_unit-sep-cov-faith}(iii),
it suffices to show assertion (i) of the proposition. To this aim,
we consider the essentially commutative diagram
$$
\xymatrix{ \Fib(u)_!(\sC(\cE)) \ar[rrr]^-{\Fib(u)_!(\sC(f))}
\ar[d]_{\Fib(u)_!(\sev^\cE)} & & &
\Fib(u)_!(\sC(\cF)) \ar[d]^{\Fib(u)_!(\sev^\cF)} \\
\Fib(u)_!(\cE) \ar[rrr]^-{\Fib(u)_!(f)} & & & \Fib(u)_!(\cF)
}$$
where $\sev^\bullet$ is the pseudo-natural equivalence described
in \eqref{subsec_eval-functor}. Then both vertical arrows of the
diagram are equivalences (lemma \ref{lem_equiv-is-preserved}(ii)),
so $\Fib(u)_!(f)$ is $i$-covering if and only if the same holds
for $\Fib(u)_!(\sC(f))$ (lemma \ref{lem_reduce-to-split-fibs}),
and by the same token, $f$ is $i$-covering if and only if the
same holds for $\sC(f)$. Thus, we may replace $f$ by $\sC(f)$,
and assume from start that $f$ is a split cartesian functor
$(\cE,\blambda^\cE)\to(\cF,\blambda^\cF)$. In this case, let
$\sc^\cE$ and $\sc^\cF$ be the strict pseudo-functors associated
with the cleavages $\blambda^\cE$ and $\blambda^\cF$; then
$\Fib(u)_!(\cE)$ is equivalent to the fibration associated with
the pseudo-functor
$$
\sd^\cE:\cC^o\to\bCat
\qquad
X\mapsto\Pscolim{(X/u\cC')^o}\sc^\cE\circ\st^o_X
$$
where $\st_X:X/u\cC'\to\cC'$ is the target functor. Hence,
$\sd^\cE_X=\cFib(\sc^\cE\circ\st^o_X)[\Sigma^{-1}_{\cE,X}]$, where
$\Sigma_{\cE,X}$ is the set of cartesian morphisms of
$\cFib(\sc^\cE\circ\st^o_X)$ (see the proof of theorem
\ref{th_bCat-cpt}). The universal pseudo-cocone
$\pi^{\cE,X}:\sc^\cE\circ\st^o_X\Rightarrow\sF_{\sd^\cE_X}$ is
the strict pseudo-natural transformation that assigns to
every object $(X',\psi:X\to uX')$ of $(X/u\cC')^o$ the
functor
$$
\pi^{\cE,X}_{(X',\psi)}:\sc^\cE_{X'}\to\sd^\cE_X
\qquad
E\mapsto((X',\psi),E).
$$
To every morphism $\phi:X\to Y$ of $\cC$ we then attach the functor
$$
\sd^\cE_\phi:\sd^\cE_Y\to\sd^\cE_X
\qquad\text{such that}\qquad
\pi^{\cE,X}*\phi^{*o}=\sF_{\sd^\cE_\phi}\odot\pi^{\cE,Y}
$$
where $\phi^*:Y/u\cC'\to X/u\cC'$ is the functor induced by
$\phi$ (see example \ref{ex_pseudo-functors}(ii)); unwinding
the definitions, we find that $\sd^\cE_\phi$ is given by the rule :
$$
((X',\psi:Y\to uX'),E)\mapsto((X',\psi\circ\phi),E)
\qquad
\text{for every $((X',\psi),E)\in\Ob(\sd^\cE_Y)$}.
$$
Likewise we describe $\Fib(u)_!(\cF)$ up to equivalence;
next, the split cartesian functor $f$ corresponds to a
strict pseudo-natural transformation
$\sc^f:\sc^\cE\Rightarrow\sc^\cF$, and under the foregoing
identifications, the cartesian functor $\Fib(u)_!(f)$
corresponds to the strict pseudo-natural transformation
$\sd^f:\sd^\cE\Rightarrow\sd^\cF$ assigning to every $X\in\Ob(\cC)$
the functor :
$$
\sd^f_X:\sd^\cE_X\to\sd^\cF_X
\qquad
((X',\psi),E)\mapsto((X',\psi),\sc^f_{X'}E).
$$
The objects of $\cFib(\sd^\cE)$ are then the data $(X,X',\psi,E)$
with $X\in\Ob(\cC)$ and $((X',\psi),E)\in\Ob(\sd^\cE_X)$, and
likewise for the objects of $\cFib(\sd^\cF)$. After this
preparation, suppose that $f$ is $0$-covering, and let us
check that the same holds for $\Fib(u)_!(f)$, where the latter
is naturally identified with $\cFib(\sd^f)$. Thus, consider any
object $(X,X',\psi:X\to uX',F)$ of  $\cFib(\sd^\cF)$; according
to remark \ref{rem_realign}, we need to find a covering family
$(\phi_\lambda:Y_\lambda\to X~|~\lambda\in\Lambda)$ and for every
$\lambda\in\Lambda$ an object
$\underline E_\lambda:=((Y'_\lambda,\psi_\lambda),E)\in\sd^\cE_{Y_\lambda}$
with an isomorphism
$$
((Y'_\lambda,\psi_\lambda),\sc^f_{Y'_\lambda}E_\lambda)=
\sd^f_{Y_\lambda}((Y'_\lambda,\psi_\lambda),E_\lambda)\isom
\sd^\cF_{\phi_\lambda}((X',\psi),F)=((X',\psi\circ\phi_\lambda),F)
\qquad
\text{in $\sd^\cF_{Y_\lambda}$}.
$$
But since $f$ is $0$-covering, there exist a covering
$(\tau_\lambda:Y'_\lambda\to X'~|~\lambda\in\Lambda)$ for the
topology $J'$ and an object $E_\lambda\in\Ob(\sc^\cE_{Y'_\lambda})$
with an isomorphism
$\omega_\lambda:\sc^f_{Y'_\lambda}E_\lambda\isom\sc^\cF_{\tau_\lambda}F$
in $\sc^\cF_{Y'_\lambda}$ for every $\lambda\in\Lambda$. Thus,
let us set $Y_\lambda:=X\times_{uX'}uY'_\lambda$, and denote by
$\phi_\lambda:Y_\lambda\to X$ the induced projection, for every
$\lambda\in\Lambda$. Since $u$ is continuous, the family
$(u(\tau_\lambda)~|~\lambda\in\Lambda)$ covers $uX'$ for the
topology $J$, and therefore $(\phi_\lambda~|~\lambda\in\Lambda)$
covers $X$. Lastly, the sought isomorphism shall be the composition
of the isomorphism $((Y'_\lambda,\psi_\lambda),\sc^f_{Y'_\lambda}E_\lambda)
\isom((Y'_\lambda,\psi_\lambda),\sc^\cF_{\tau_\lambda}F)$ of
$\cFib(\sc^\cF\circ\st^o_{Y_\lambda})$ induced by $\omega_\lambda$,
and the cartesian morphism
$((Y'_\lambda,\psi_\lambda),\sc^\cF_{\tau_\lambda}F)\to
((X',\psi\circ\phi_\lambda),F)$ of $\cFib(\sc^\cF\circ\st^o_{Y_\lambda})$.

Next, suppose that $f$ is $1$-covering, and let
$((X'_i,\psi_i:X\to uX'_i),E_i)$ for $i=1,2$ be two objects of
$\sd^\cE_X$ with a morphism $g:((X'_1,\psi_1),\sc^f_{X'_1}E_1)
\to((X'_2,\psi_2),\sc^f_{X'_2}E_2)$ in $\sd^\cF_X$. Notice that
$(Z/u\cC')^o$ is a filtered category for every $Z\in\Ob(\cC)$
(example \ref{ex_cofiltered-comma}(i)), hence $\Sigma_{\cF,Z}$
admits a right calculus of fraction (example
\ref{ex_filter-2-colim-in-Cat}(ii)); moreover, every cartesian
morphism $(\sigma,s):((Z'',\psi'),F')\to((Z',\psi),F)$ in
$\Sigma_{\cF,Z}$ admits the factorization :
\set\begin{equation}\label{eq_six-in-Jerusalem}
(\sigma,s)=\blambda^\cF(F,\sigma)\circ(\one_{Z''},s)
\qquad
\text{in $\cFib(\sc^\cF\circ\st^o_Z)$}
\end{equation}
where $(\one_{X''},s)$ is an isomorphism of $\sc^\cF_{X''}$
(lemma \ref{lem_distinguished-cleavage}(i)), so $g$ can be
represented by a pair
$$
((X'_1,\psi_1),\sc^f_{X'_1}E_1)
\xleftarrow{\blambda^\cF(\sc^f_{X'_1}E_1,\rho)}
((X''_1,\psi'_1),\sc^\cF_\rho\circ\sc^f_{X'_1}E_1)\xrightarrow{g'}
((X'_2,\psi_2),\sc^f_{X'_2}E_2)
$$
where $\rho:X''_1\to X'_1$ is a morphism of $\cC'$ and
$\psi'_1:X\to uX''_1$ is a morphism of $\cC$ such that
$u(\rho)\circ\psi'_1=\psi_1$. Also, $g'$ is a
morphism in $\cFib(\sc^\cF\circ\st^o_X)$. Then in turn,
$g'$ is a pair $(X/\mu,t)$, where $X/\mu$ is a morphism
$(X''_1,\psi'_1)\to(X'_2,\psi_2)$ in $X/u\cC'$, {\em i.e.}
a morphism $\mu:X''_1\to X'_2$ in $\cC'$ with
$u(\mu)\circ\psi_1=\psi_2$, and $g'$ is a morphism in
$\sc^\cF_{X''_1}$
$$
g':\sc^f_{X''_1}\circ\sc^\cE_\rho E_1=\sc^\cF_\rho\circ\sc^f_{X'_1}E_1
\to\sc^\cF_\mu\circ\sc^f_{X'_2}E_2=\sc^f_{X''_1}\circ\sc^\cE_\mu E_2.
$$
According to remark \ref{rem_realign}, we need to exhibit
a covering family
$(\phi_\lambda:Y_\lambda\to X~|~\lambda\in\Lambda)$, and
for every $\lambda\in\Lambda$ a morphism
$g_\lambda:((X'_1,\psi_1\circ\phi_\lambda),E_1)\to
((X'_2,\psi_2\circ\phi_\lambda),E_2)$ in $\sd^\cE_{Y_\lambda}$
such that $\sd^f_{X'_1}(g_\lambda)=\sd^\cF_{\phi_\lambda}(g)$.
Now, since $f$ is $1$-covering, there exist a covering
family $(\tau_\lambda:Y''_\lambda\to X''_1~|~\lambda\in\Lambda)$
and for every $\lambda\in\Lambda$ a morphism
$$
g'_\lambda:\sc^\cE_{\tau_\lambda}\circ\sc^\cE_\rho E_1\to
\sc^\cE_{\tau_\lambda}\circ\sc^\cE_\mu E_2
\qquad\text{in $\sc^\cE_{Y''_\lambda}$ such that}\qquad
\sc^f_{Y''_\lambda}(g'_\lambda)=\sc^\cF_{\tau_\lambda}(g').
$$
We set $Y_\lambda:=X\times_{uX''_1}uY''_\lambda$, and we
let $\phi_\lambda:Y_\lambda\to X$ and
$\psi''_\lambda:Y_\lambda\to uY''_\lambda$ be the induced
projections for every $\lambda\in\Lambda$; arguing as
in the foregoing it is easily seen that the family
$(g_\lambda~|~\lambda\in\Lambda)$ covers $X$, and we
define $g_\lambda$ as the morphism of $\sd^\cE_{Y_\lambda}$
represented by the pair of morphisms :
$$
((X'_\lambda,\psi_1\circ\phi_\lambda),E_1)
\xleftarrow{\blambda^\cE(E_1,\rho\circ\tau_\lambda)}
((Y''_\lambda,\psi''_\lambda),\sc^\cE_{\rho\circ\tau_\lambda}E_1)
\xrightarrow{(Y_\lambda/\mu\circ\tau_\lambda,g'_\lambda)}
((X'_2,\psi_2\circ\phi_\lambda),E_2).
$$
In view of remark \ref{rem_conditions-for-split}(ii) we
see that $\sd^f_{Y_\lambda}(g'_\lambda)$ is represented by the pair
$$
((X'_\lambda,\psi_1\phi_\lambda),\sc^f_{X'_1}E_1)
\xleftarrow{\blambda^\cF(\sc^f_{X'_1}E_1,\rho\tau_\lambda)}
((Y''_\lambda,\psi''_\lambda),\sc^f_{Y''_\lambda}\sc^\cE_{\rho\tau_\lambda}E_1)
\xrightarrow{(Y_\lambda/\mu\tau_\lambda,\sc^\cF_{\tau_\lambda}(g'))}
((X'_2,\psi_2\phi_\lambda),\sc^f_{X'_2}E_2)
$$
and by remark \ref{rem_conditions-for-split}(i) we have :
$$
\blambda^\cF(\sc^f_{X'_1}E_1,\rho\tau_\lambda)=
\blambda^\cF(\sc^f_{X'_1}E_1,\rho)\circ
\blambda^\cF(\sc^\cF_\rho\sc^f_{X'_1}E_1,\tau_\lambda).
$$
On the other hand,
$(Y_\lambda/\mu\circ\tau_\lambda,\sc^\cF_{\tau_\lambda}(g'))$ is
the composition :
$$
((Y''_\lambda\!,\!\psi''_\lambda),\sc^f_{Y''_\lambda}\sc^\cE_{\rho\tau_\lambda}\!E_1)
\xrightarrow{\blambda^\cF(\sc^\cF_\rho\sc^f_{X'_1}E_1,\tau_\lambda)}
((X''_\lambda,\psi'_1\phi_\lambda),\sc^\cF_\rho\sc^f_{X'_1}E_1)
\xrightarrow{\sd^\cF_{\phi_\lambda}(X/\mu,g')}
((X'_2,\psi_2\phi_\lambda),\sc^f_{X'_2}E_1)
$$
so $\sd^f_{Y_\lambda}(g'_\lambda)$ is also represented by the pair
$$
((X'_\lambda,\psi_1\phi_\lambda),\sc^f_{X'_1}E_1)
\xleftarrow{\blambda^\cF(\sc^f_{X'_1}E_1,\rho)}
((X''_\lambda,\psi'_1\phi_\lambda),\sc^\cF_\rho\sc^f_{X'_1}E_1)
\xrightarrow{\sd^\cF_{\phi_\lambda}(X/\mu,g')}
((X'_2,\psi_2\phi_\lambda),\sc^f_{X'_2}E_1)
$$
which represents $\sd^\cF_{\phi_\lambda}(g)$ as well, as required.

Lastly, suppose that $f$ is $2$-covering, and consider objects
$((X'_j,\psi_j),E_j)$ of $\sd^\cE_X$ for $j=1,2$, and two morphisms
$g_1,g_2:((X'_1,\psi_1),E_1)\to((X'_2,\psi_2),E_2)$ in $\sd^\cE_X$
such that
\set\begin{equation}\label{eq_campus-punk}
\sd^f_X(g_1)=\sd^f_X(g_2).
\end{equation}
We need to exhibit a covering family
$(\phi_\lambda:Y_\lambda\to X~|~\lambda\in\Lambda)$ such that
$\sd^\cE_{\phi_\lambda}(g_1)=\sd^\cE_{\phi_\lambda}(g_2)$ for every
$\lambda\in\Lambda$. We reduce easily to the case where $g_1$ and
$g_2$ are two classes of morphisms of $\cFib(\sc^\cE\circ\st^o_X)$.
In this case, \eqref{eq_campus-punk} means that there exists a
cartesian morphism
$h:((X''_1,\psi'_1),F)\to((X'_1,\psi_1),\sc^f_{X'_1}E_1)$ in
$\cFib(\sc^\cF\circ\st^o_X)$ such that
$\cFib(\sc^f\circ\st^o_X)(g_1)\circ h=
\cFib(\sc^f\circ\st^o_X)(g_2)\circ h$. However, we have a
factorization
$h=\blambda^\cF(\sc^f_{X'_1}E_1,\mu)\circ(\one_{X''_1},h')$
as in \eqref{eq_six-in-Jerusalem}, where $(\one_{X''_1},h')$
is an isomorphism of $\sc^\cF_{X''_1}$, so we may assume that
$h=\blambda^\cF(\sc^f_{X'_1}E_1,\mu)=
\cFib(\sc^f\circ\st^o_X)(\blambda^\cE(E_1,\mu))$. Then we may
as well replace $g_i$ by $g_i\circ\blambda^\cE(E_1,\mu)$ for
$i=1,2$, and assume from start that
$$
\cFib(\sc^f\circ\st^o_X)(g_1)=\cFib(\sc^f\circ\st^o_X)(g_2).
$$
Especially, there exist a morphism $\rho:X'_1\to X'_2$ in
$\cC'$ such that $u(\rho)\circ\psi_1=\psi_2$, and morphisms
$g'_i:E_1\to\sc^\cE_\rho E_2$ such that $g_i=(\rho,g'_i)$ for
$i=1,2$, and $\sc^f_{X'_1}(g'_1)=\sc^f_{X'_1}(g'_2)$. Then,
since $f$ is $2$-covering, we find a covering family
$(\tau_\lambda:Y'_\lambda\to X'_1~|~\lambda\in\Lambda)$ such
that
\set\begin{equation}\label{eq_amelie-list-to-muzak}
\sc^\cE_{\tau_\lambda}(g'_1)=\sc^\cE_{\tau_\lambda}(g'_2)
\qquad
\text{for every $\lambda\in\Lambda$}.
\end{equation}
We set again $Y_\lambda:=X\times_{uX'_1}uY'_\lambda$, and we let
$\phi_\lambda:Y_\lambda\to X$ and
$\psi_\lambda:Y_\lambda\to Y'_\lambda$ be the induced projections,
for every $\lambda\in\Lambda$. We claim that the resulting
covering $(\phi_\lambda~|~\lambda\in\Lambda)$ of $X$ will do.
Indeed, $\sd^\cE_{\phi_\lambda}(g_i):((X'_1,\psi_1\circ\phi_\lambda),E_1)
\to((X'_2,\psi_2\circ\phi_\lambda),E_2)$ is represented by the
morphism $(Y_\lambda/\rho,g'_i)$ of $\cFib(\sc^\cE\circ\st^o_{Y_\lambda})$
for $i=1,2$, and \eqref{eq_amelie-list-to-muzak} implies easily that
$$
(Y_\lambda/\rho,g'_1)\circ\blambda^\cE(E_1,\tau_\lambda)=
(Y_\lambda/\rho,g'_2)\circ\blambda^\cE(E_1,\tau_\lambda)
$$
where $\blambda^\cE(E_1,\tau_\lambda):
((Y'_\lambda,\psi_\lambda),\sc^\cE_{\tau_\lambda}E_1)\to
((X'_1,\psi_1\circ\phi_\lambda),E_1)$ is a cartesian morphism
of the fibration $\cFib(\sc^\cE\circ\st^o_{Y_\lambda})$, whence
the contention.
\end{pfclaim}

We consider now the essentially commutative diagram
$$
\xymatrix{ \Fib(u)_!(\cE) \ar[rr]^-{\Fib(u)_!(f)}
\ar[d]_{\Fib(u)_!(\eta_\cE)} & &
\Fib(u)_!(\cF) \ar[d]^{\Fib(u)_!(\eta_\cF)} \\
\Fib(u)_!(\cE^a) \ar[rr]^-{\Fib(u)_!(f^a)} & &
\Fib(u)_!(\cF^a)
}$$
where $\eta_\cE:\cE\to\cE^a$ and $\eta_\cF:\cF\to\cF^a$ are the
natural morphisms, and notice that $\Fib(u)_!(\eta_\cE)^a$ and
$\Fib(u)_!(\eta_\cF)^a$ are equivalences, by virtue of proposition
\ref{prop_weak-morph-of-sites}. Then $\Fib(u)_!(f)^a$ is
$i$-covering if and only if the same holds for $\Fib(u)_!(f^a)^a$
(lemma \ref{lem_reduce-to-split-fibs}), and finally
$\Fib(u)_!(f)$ is $i$-covering if and only if the same holds for
$\Fib(u)_!(f^a)$ (proposition \ref{prop_unit-sep-cov-faith}(iii)).
Summing up, we see that, in order to prove the proposition,
we may assume that $\cE$ and $\cF$ are stacks on $C'$, and
moreover, in this case it suffices to show assertion (ii).

Now, by lemma \ref{lem_cont-funct-site}(ii) we have the
essentially commutative diagram :
$$
\xymatrix{ \cC' \ar[r]^-u \ar[d]_{h^a_{C'}} &
\cC \ar[d]^{h^a_C} \\
C'^\sim_\sU \ar[r]^-{\tilde u{}^*} & C^\sim_\sU
}$$
whose vertical arrows are the Yoneda embeddings. Especially,
for every universe $\sV$ containing $\sU$ and such that
$\cC^\wedge_\sU$ is $\sV$-small, $\sV\tdu\St(h^a_{C'})_*$ is
a $2$-equivalence from the $2$-category of stacks on
$\Can(C'^\sim_\sU)$ to that of stacks on $C$ (theorem
\ref{th_2-equiv-for-cats-of-stacks}), and therefore
there exists a cartesian functor $f':\cE'\to\cF'$ of
$\sV$-stacks on $\Can(C'^\sim_\sU)$ with an essentially
commutative diagram :
$$
\cD
\quad :\quad
{\diagram \cE \ar[rrr]^-f \ar[d] & & & \cF \ar[d] \\
\sV\tdu\St(h^a_{C'})_*(\cE') \ar[rrr]^-{\sV\tdu\St(h^a_{C'})_*(f')}
& & & \sV\tdu\St(h^a_{C'})_*(\cF')
\enddiagram}$$
whose vertical arrows are equivalences. We notice :

\begin{claim}\label{cl_more-general-than-need}
Let $C:=(\cC,J)$ be a small site, $h^a_C:\cC\to C^\sim$
the Yoneda embedding, $\phi:\cA\to\cB$ a morphism of
$\sV\tdu\Fib(C^\sim)$, and $j\in\{0,1,2\}$. Then $\phi$
is $j$-covering for the canonical topology of $C^\sim$
if and only if $\sV\tdu\Fib(h^a_C)^*(\phi)$ is $j$-covering
for the topology $J$.
\end{claim}
\begin{pfclaim} Suppose that $\phi$ is $j$-covering;
since $h^a_C$ is cocontinuous for the topologies $J$
(theorem \ref{th_canon-topos}(iii)) and $\Can_{C^\sim}$,
proposition \ref{prop_cocont-and-stacks}(i) says that
$\sV\tdu\Fib(h^a_C)^*(\phi)$ is $j$-covering. The converse
follows from lemma \ref{lem_kayak-days}.
\end{pfclaim}

Since $f$ is $i$-covering, the same holds for
$\sV\tdu\St(h^a_{C'})_*(f')$, by lemma \ref{lem_reduce-to-split-fibs},
and then also for $f'$, by claim \ref{cl_more-general-than-need}.
Moreover, considering the essentially commutative diagram
$\sV\tdu\St(u)^*(\cD)$ and applying again lemma
\ref{lem_reduce-to-split-fibs}, we are reduced to checking that
$\sV\tdu\St(u)^*\circ\sV\tdu\St(h^a_{C'})_*(f')$
is $i$-covering. To ease notation, set
$S:=\sV\tdu\St(h^a_C)_*\circ\sV\tdu\St(h^a_C)^*\circ
\sV\tdu\St(u)^*\circ\sV\tdu\St(h^a_{C'})_*$;
recalling that $\sV\tdu\St(h^a_C)^*$ is a pseudo-inverse for
the $2$-equivalence $\sV\tdu\St(h^a_C)_*$, we get an essentially
commutative diagram
$$
\xymatrix@C+10pt{ \sV\tdu\St(u)^*\circ\sV\tdu\St(h^a_{C'})_*(\cE')
\ar[rrr]^-{\sV\tdu\St(u)^*\circ\sV\tdu\St(h^a_{C'})_*(f')} \ar[d]
& & & \sV\tdu\St(u)^*\circ\sV\tdu\St(h^a_{C'})_*(\cF') \ar[d] \\
S(\cE) \ar[rrr]^-{S(f')} & & & S(\cF)
}$$
whose vertical arrows are equivalences; invoking once again
lemma \ref{lem_reduce-to-split-fibs}, we are then reduced
to checking that $S(f)$ is $i$-covering. However, we have
pseudo-natural equivalences :
$$
S\isom\sV\tdu\St(h^a_C)_*\circ\sV\tdu\St(\tilde u^*)^*\circ
\sV\tdu\St(h^a_{C'})^*\circ\sV\tdu\St(h^a_{C'})_*\isom
\sV\tdu\St(h^a_C)_*\circ\sV\tdu\St(\tilde u^*)^*
$$
which, after applying once more lemma \ref{lem_reduce-to-split-fibs},
further reduces to checking that the functor
$\sV\tdu\St(h^a_C)_*\circ\sV\tdu\St(\tilde u^*)^*(f')$ is
$i$-covering. The latter follows from claims
\ref{cl_ok-for-lex-sites} and \ref{cl_more-general-than-need}.
\end{proof}

\sset\subsubsection{}\label{subsec_stupeur}
Let $C:=(\cC,J)$ and $C':=(\cC',J')$ be two small sites,
and $u:\cC\to\cC'$ a cocontinuous functor; by proposition
\ref{prop_cocont-and-stacks}(ii), the pseudo-functor
$\Fib(u)_*$ restricts to a pseudo-functor
$$
\St(\breve u)_*:\Stack(C)\to\Stack(C')
$$
and it is easily seen that the latter admits a left $2$-adjoint
$$
\St(\breve u)^*:\Stack(C')\xrightarrow{i}\Fib(\cC')
\xrightarrow{\Fib(u)^*}\Fib(\cC)\xrightarrow{(-)^a_C}\Stack(C)
$$
with $i$ the inclusion pseudo-functor. On the other hand, by
remark \ref{rem_choose-two-univs}(iii,iv), the functor $u$ induces
a morphism of sites $\breve u{}^*:\Can(C^\sim)\to\Can(C'^\sim)$,
whence a $2$-adjoint pair of pseudo-functors
$(\St(\breve u{}^*)^*,\St(\breve u{}^*)_*)$. Define also
$h^a_C:\cC\to C^\sim,h^a_{C'}:\cC'\to C'^\sim$ as in remark
\ref{rem_rep-and-sheafify}(iii).

\begin{proposition}\label{prop_breve-for-stacks}
{\em (i)}\ \
In the situation of \eqref{subsec_stupeur}, we have
pseudo-commutative diagrams :
$$
\xymatrix@C+6pt{
\Stack(\Can\ C'^\sim) \ar[r]^-{\St(\breve u{}^*)^*}
\ar[d]_{\St(h^a_{C'})_*} &
\Stack(\Can\ C^\sim) \ar[d]^{\St(h^a_C)_*} &
\Stack(\Can\ C^\sim) \ar[r]^-{\St(\breve u{}^*)_*}
\ar[d]_{\St(h^a_C)_*} &
\Stack(\Can\ C'^\sim) \ar[d]^{\St(h^a_{C'})_*} \\
\Stack(C') \ar[r]^-{\St(\breve u)^*} & \Stack(C) &
\Stack(C) \ar[r]^-{\St(\breve u)_*} & \Stack(C').
}$$

{\em (ii)}\ \
If $u$ is also a weak morphism of sites $C'\to C$, we have a
pseudo-natural equivalence
$$
\St(\breve u)^*\isom\St(u)_*.
$$

{\em (iii)}\ \
If $u$ admits a right adjoint $v:\cC'\to\cC$ that is a weak
morphism of sites $C\to C'$, we have pseudo-natural equivalences :
$$
\St(v)_*\isom\St(\breve u)_*
\qquad
\St(v)^*\isom\St(\breve u)^*.
$$
\end{proposition}
\begin{proof}(i): Since $\St(h^a_C)_*$ and $\St(h^a_{C'})_*$ are
$2$-equivalences (theorem \ref{th_2-equiv-for-cats-of-stacks}(i,ii)),
it suffices to show the pseudo-commutativity of the left square
diagram. To this aim, consider the sites $C^\wedge$ and $C'^\wedge$
on the categories $\cC^\wedge$ and $\cC'^\wedge$ defined in remark
\ref{rem_topol-on-presheaves}(ii), and notice that since
$u^\wedge:\cC^\wedge\to\cC'^\wedge$ commutes with all limits
and all colimits, it induces a morphism of topoi
$C'^\wedge\to C^\wedge$ (proposition \ref{prop_half-dominates}(i)),
whence a morphism of sites $u^\wedge:C'^\wedge)\to C^\wedge)$
(remark \ref{rem_choose-two-univs}(iii)); combining with lemma
\ref{lem_breve} and theorem \ref{th_canon-topos}(ii) we
deduce an essentially commutative diagram of sites:
$$
\xymatrix{ \Can(C'^\sim) \ar[r]^-{\breve u{}^*} \ar[d]_{(-)_{C'}^a} &
\Can(C^\sim) \ar[d]^{(-)^a_C} \\
C'^\wedge \ar[r]^-{u^\wedge} & C^\wedge.
}$$
Notice as well that
$\St((-)_C^a)_*:\Stack(\Can(C^\sim))\to\Stack(C^\wedge)$ is a
$2$-equivalence (proposition \ref{prop_equiv-and-2-equiv}(ii)
and theorem \ref{th_canon-topos}(iii)), and likewise for
$\St((-)_{C'}^a)_*$. There follows a pseudo-commutative diagram
(details left to the reader) :
$$
\xymatrix@C+20pt{ \Stack(\Can(C'^\sim))
\ar[r]^-{\St(\breve u{}^*)^*} \ar[d]_{\St((-)_{C'}^a)_*} &
\Stack(\Can(C^\sim)) \ar[d]^{\St((-)^a_C)_*} \\
\Stack(C'^\wedge) \ar[r]^-{\St(u^\wedge)^*} & \Stack(C^\wedge).
}$$
Thus, we are reduced to checking that the following diagram
pseudo-commutes :
$$
\xymatrix@C+40pt{ \Fib(\cC'^\wedge) \ar[r]^-{\Fib(u^\wedge)_!}
\ar[d]_{\Fib(h_{\cC'})^*} & \Fib(\cC^\wedge) \ar[r]^-{(-)^a_{C^\wedge}}
\ar[d]_{\Fib(h_\cC)^*} & \Stack(C^\wedge) \ar[d]^{\St(h_\cC)_*} \\
\Fib(\cC') \ar[r]^-{\Fib(u)^*} &
\Fib(\cC) \ar[r]^-{(-)^a_C} & \Stack(C).
}$$
However, the pseudo-commutativity of the left square
subdiagram follows from proposition \ref{prop_same-vein-Fib}.
For the right subdiagram, it suffices to observe that
$h_\cC$ is cocontinuous for the sites $C$ and $C^\wedge$
(theorem \ref{th_canon-topos}(i)) and invoke proposition
\ref{prop_cocont-and-stacks}(i).

(ii) follows by inspecting the definitions, and the first
pseudo-natural equivalence of (iii) follows from remark
\ref{rem_explicit-Fib_*}(vi); then the second pseudo-natural
equivalence of (iii) follows from the first one and from
remark \ref{rem_opposing-thumb}(iii).
\end{proof}

\sset\subsubsection{}
\label{subsec_weak-morph-of-sites-and-cocont}
Let $C:=(\cC,J)$ and $C':=(\cC',J')$ be two $\sU$-sites,
and $u:C'\to C$ a weak morphism of sites. We have a
commutative diagram of $2$-categories :
$$
\cD_u \qquad : \qquad
{\diagram \Stack(C') \ar[r]^-{i_{C'}} \ar[d]_{\St(u)_*} &
\Fib(\cC') \ar[d]^{\Fib(u)^*} \\
\Stack(C) \ar[r]^-{i_C} & \Fib(\cC)
\enddiagram}
\qquad\qquad$$
where the inclusion pseudo-functors $i_C$ and $i_{C'}$ admit
left $2$-adjoints $(-)^a_C:\Fib(\cC)\to\Stack(C)$ and
$(-)^a_{C'}:\Fib(\cC')\to\Stack(C')$ (corollary
\ref{cor_stackificate-U-sites}). After fixing an adjunction
for these two pairs of $2$-adjoint pseudo-functors, we may then
regard $\cD_u$ as a square of weak links as in
\eqref{subsec_variant-of-Upsilon}, oriented by the identity
pseudo-natural transformation. If we similarly fix adjunctions
also for the pairs $(\Fib(u)_!,\Fib(u)^*)$ and $(\St(u)^*,\St(u)_*)$,
the diagram $\cD_u$ becomes even an oriented square of links.

\begin{corollary}\label{cor_weak-morph-of-sites-and-cocont}
In the situation of \eqref{subsec_weak-morph-of-sites-and-cocont},
suppose that the functor $u:\cC\to\cC'$ is also cocontinuous
for the topologies $J$ and $J'$. Then the base change
transformation
$$
\Upsilon(\cD_u):(-)^a_C\circ\Fib(u)^*\to\St(u)_*\circ(-)^a_{C'}
$$
is a pseudo-natural equivalence.
\end{corollary}
\begin{proof} Let $\sV$ be a universe with $\sU\subset\sV$,
and such that $\cC$ and $\cC'$ are $\sV$-small; we get a
commutative diagram of $2$-categories :
$$
\xymatrix{ \sV\tdu\Stack(C') \ar[rrr]^-{\sV\tdu i_{C'}}
\ar[ddd]_{\sV\tdu\St(u)_*} & & & \sV\tdu\Fib(\cC')
\ar[ddd]^{\sV\tdu\Fib(u)^*} \\
& \Stack(C') \ar[r]^-{i_{C'}} \ar[d]_{\St(u)_*} \ar[lu]_{j_{C'}} &
\Fib(\cC') \ar[d]^{\Fib(u)^*} \ar[ru]^{j_{\cC'}} \\
& \Stack(C) \ar[r]^-{i_C} \ar[ld]_{j_C} & \Fib(\cC) \ar[rd]^{j_\cC} \\
\sV\tdu\Stack(C) \ar[rrr]^-{\sV\tdu i_C} & & & \sV\tdu\Fib(\cC)
}$$
whose horizontal arrows are the inclusion pseudo-functors;
since the latter all admit left $2$-adjoints, we may regard
this diagram as an oriented cubical diagram of weak links as
in \ref{subsec_transfer-base-change}, whose orientations
are given by identities. Especially, the lower and upper
trapezoidal subdiagrams are oriented respectively by
$\one_{j_\cC\circ i_C}$ and $\one_{j_{\cC'}\circ i_{C'}}$, and we remark :

\begin{claim}\label{cl_fila-via}
$\Upsilon(\one_{j_\cC\circ i_C})$ and
$\Upsilon(\one_{j_{\cC'}\circ i_{C'}})$ are pseudo-natural
equivalences.
\end{claim}
\begin{pfclaim} We show the assertion for
$\Upsilon(\one_{j_\cC\circ i_C})$ : the same argument shall apply
to $\Upsilon(\one_{j_{\cC'}\circ i_{C'}})$ as well. Let $(\eta,\eps)$
(resp. $(\sV\tdu\eta,\sV\tdu\eps)$) be the unit and counit of
a $2$-adjunction for the $2$-adjoint pair $((-)^a_C,i_C)$ (resp.
for the $2$-adjoint pair $(\sV\tdu(-)^a_C,\sV\tdu i_C)$); since
$\sV\tdu i_C$ is fully faithful, $\sV\tdu\eps$ is a pseudo-natural
equivalence, so it suffices to check that the same holds for
$(\sV\tdu(-)^a_C\circ j_\cC)*\eta$. But in view of the pseudo-natural
equivalence $\sV\tdu(-)^a_C\circ j_\cC\isom j_C\circ(-)^a_C$, we
are further reduced to checking that $(-)_C^a*\eta$ is a
pseudo-natural equivalence; the latter follows from the
triangular identities for the pair $(\eta,\eps)$, since $\eps$
is a pseudo-natural equivalence.
\end{pfclaim}

Denote by $\sV\tdu\cD_u$ the ``front face'' of the foregoing
diagram (with its identity orientation); from claim
\ref{cl_fila-via} and remark \ref{rem_transit-base-change}(i)
we see that if $\Upsilon(\sV\tdu\cD_u)$ is a pseudo-natural
equivalence, the same holds for $\Upsilon(\cD_u)$. Thus, we
may replace $\sU$ by $\sV$, and assume from start that $\cC$
and $\cC'$ are small. Let now $\cF$ be any fibration on $\cC'$;
the functor $\Upsilon(\cD_u)_\cF:(\Fib(u)^*\cF)^a\to\St(u)_*(\cF^a)$
is obtained explicitly as follows. First, let
$\eta_\cF:\cF\to\cF^a$ be the unit of the chosen $2$-adjunction
for the pair $((-)^a_{C'},i_{C'})$. Then $\Upsilon(\cD_u)_\cF$ is
the composition of the induced functor
$$
(\Fib(u)^*\eta_\cF)^a:(\Fib(u)^*\cF)^a\to(\St(u)_*(\cF^a))^a
$$
with the counit of adjunction
$\eps_{\cF^a}:(\St(u)_*(\cF^a))^a\to\St(u)_*(\cF^a)$. We need
to check that $\Upsilon(\cD_u)_\cF$ is an equivalence; however,
$\eps_{\cF^a}$ is an equivalence, so it suffices to check
that the same holds for $(\Fib(u)^*\eta_\cF)^a$. By propositions
\ref{prop_faith-separ-cover} and \ref{prop_unit-sep-cov-faith}(iii)
we are reduced to showing that
$\Fib(u)^*(\eta_\cF):\Fib(u)^*\cF\to\St(u)_*(\cF^a)$ is
$i$-covering for $i=0,1,2$. Since $u$ is cocontinuous,
the latter assertion follows from propositions
\ref{prop_cocont-and-stacks} and \ref{prop_unit-sep-cov-faith}(i).
\end{proof}

\subsection{Sheaves of categories}
We wish next to elucidate the relationship between stacks and
the related notion of sheaf of categories on a given site. We
begin by recalling the following :

\begin{definition}\label{def_sheaves-with-other-values}
Let $C:=(\cC,J)$ be any site, and $\cA$ any other category.

(i)\ \
A {\em presheaf on $\cC$ with values in $\cA$} is any object
of the category (notation of \eqref{subsec_Godem-prod})
$$
(\cC,\cA)^\wedge:=\bFun(\cC^o,\cA).
$$

(ii)\ \
Let $\sV$ be a universe such that $\cA$ has $\sV$-small
$\Hom$-sets; following \cite[Ch.0, \S3.1]{EGAI}, we say
that such a presheaf $\cF$ is a {\em sheaf on $C$ with
values in $\cA$} if for every $T\in\Ob(\cA)$, the rule :
$$
U\mapsto\cF_T(U):=\Hom_\cA(T,\cF(U))
\qquad
\text{for every $U\in\Ob(\cC)$}
$$
defines a $\sV$-sheaf (of sets) $\cF_T$ on $C$. This condition
is obviously independent of the choice of $\sV$. We denote by
$$
(C,\cA)^\sim
$$
the full subcategory of $(\cC,\cA)^\wedge$ whose objects
are the sheaves on $C$ with values in $\cA$.
\end{definition}

\begin{remark}\label{rem_sheaves-with-values-in-A}
(i)\ \ 
Let $C:=(\cC,J)$ be a site and $\cF$ a presheaf on $\cC$ with
values in a category $\cA$. According to
\eqref{subsec_interpret-descent}, the presheaf $\cF$ is a sheaf
on $C$ with values in $\cA$ if and only if, for every
$T\in\Ob(\cA)$, every $U\in\Ob(\cC)$, and every sieve
$\cS\in J(U)$ the natural map
$$
\cF_T(U)\to\lim_{\cS^o}\cF_T\circ\ss^o_\cS
$$
is bijective (where the limit is taken in a sufficiently large
universe $\sV$). The latter in turns is equivalent to the following.
For every $U$ and $\cS$ as in the foregoing, the induced cone :
$$
(\cF(f):\cF(U)\to\cF(U')~|~(f:U'\to U)\in\Ob(\cS))
$$
is universal in $\cA$. Especially, a sheaf on $C$ with values in
$\Set$ is just a usual sheaf (of sets).

In case $\cC$ is small, $\cA$ is complete, and all the fibre
products in $\cC$ are representable, arguing similarly we see
that $\cF$ is a sheaf on $C$ with values in $\cA$ if and only
if, for every $U\in\Ob(\cC)$ and every small family
$(U_i\to U~|~i\in I)$ of objects of $\cC/U$ that generate a sieve
covering $U$, the restriction morphisms $\cF(U)\to\cF(U_i)$ induce
an isomorphism
$$
\cF(U)\isom\Equal\Bigl(\prod_{i\in I}
\xymatrix{\cF(U_i)\ar@<-.5ex>[r] \ar@<.5ex>[r] &} \!\!\!\!\!\!\!
\prod_{(i,j)\in I\times I}\!\!\!\cF(U_i\times_UU_j)\Bigr)
\qquad
\text{in $\cA$}.
$$

(ii)\ \
Let $\cA'$ be another category, and $F:\cA\to\cA'$ a functor
that commutes with limits; it follows from (i) that
$(\cC,F)^\wedge:=\bFun(\cC^o,F):(\cC,\cA)^\wedge\to(\cC,\cA')^\wedge$
restricts to a functor
$$
(C,F)^\sim:(C,\cA)^\sim\to(C,\cA')^\sim.
$$

(iii)\ \
Suppose that $\cA$ is complete and $\sV$-small for a universe
$\sV$ such that $\sU\subset\sV$. Consider a functor
$$
\cF_\bullet:\Lambda\to(C,\cA)^\sim
\qquad
\lambda\mapsto\cF_\lambda
$$
from a small category $\Lambda$. Recall that the limit $L$ of
$\cF_\bullet$ in $(\cC,\cA)^\wedge$ is computed argumentwise,
{\em i.e.} $L(U)$ represents the limit of the induced functor
$$
\cF_\bullet(U):\Lambda\to\cA
\qquad
\lambda\mapsto\cF_\lambda(U)
$$
for every $U\in\Ob(\cC)$ (corollary
\ref{lem_presheaf-in-a-cat}(ii)). We claim that $L$ is
a sheaf on $C$ with values in $\cA$, and therefore it
represents the limit of $\cF_\bullet$ in $(C,\cA)^\sim$;
especially, the latter category is complete. Indeed,
for any $T\in\Ob(\cA)$ we have natural identifications
$$
\Hom_\cA(T,L(U))\isom
\lim_{\lambda\in\Ob(\Lambda)}\Hom_\cA(T,\cF_\lambda(U))
=\lim_{\lambda\in\Ob(\Lambda)}\cF_{\lambda,T}(U)\isom
\Bigl(\lim_{\lambda\in\Ob(\Lambda)}\cF_{\lambda,T}\Bigr)(U)
$$
and the limit of the functor $\cF_{\bullet,T}$ from $\Lambda$
to the category of $\sV$-presheaves (of sets) on $\cC$ is a
$\sV$-sheaf, since by assumption every $\cF_{\lambda,T}$ is
a $\sV$-sheaf (remark \ref{rem_rep-and-sheafify}(i)).
This shows that $L_T$ is a $\sV$-sheaf on $C$ for every
$T\in\Ob(\cA)$, as required.

(iv)\ \
Let $C':=(\cC',J')$ be another site, and $u:\cC\to\cC'$ a
continuous functor for the topologies $J$ and $J'$. Then
$u$ induces a functor
$$
(u,\cA)^\wedge:=\bFun(u^o,\cA):(\cC',\cA)^\wedge\to(\cC,\cA)^\wedge
$$
and notice that for every presheaf $\cF'$ on $\cC'$ with values
in $\cA$, and every $T\in\Ob(\cA)$ we have
\set\begin{equation}\label{eq_restrict-sheaf}
u^\wedge(\cF'_T)=((u,\cA)^\wedge\cF')_T.
\end{equation}
It follows that $(u,\cA)^\wedge$ restricts to a functor
$$
(\tilde u,\cA)_*:(C',\cA)^\sim\to(C,\cA)^\sim.
$$
Likewise, every natural transformation $\alpha:u\Rightarrow v$
between such continuous functors $u,v:\cC\to\cC'$ induces a
natural transformation
$$
(\alpha,\cA)^\wedge:=\bFun(\alpha^o,\cA):
(v,\cA)^\wedge\Rightarrow(u,\cA)^\wedge
$$
which yields by restriction a natural transformation
$$
(\tilde \alpha,\cA)_*:(v,\cA)^\sim_*\Rightarrow(u,\cA)^\sim_*.
$$

(v)\ \
For every $U\in\Ob(\cC)$, the topology $J$ induces a topology $J_U$
on $\cC/U$ (see \eqref{subsec_topol-on-C-over-X}), and we let
$C/U:=(\cC/U,J_U)$. The source functor $\ss_U:\cC/U\to\cC$ is
continuous for the topologies $J$ and $J_U$, hence (iv) yields
a well defined functor
$$
(\tilde\ss_U,\cA)_*:(C,\cC)^\sim\to(C/U,\cC)^\sim
\qquad
\cF\mapsto\cF_{|U}.
$$
Likewise, if $g:U\to V$ is any morphism of $\cC$, the
corresponding functor $g_*:\cC/U\to\cC/V$ is continuous for
the topologies $J_U$ and $J_V$, so it restricts to a functor
$$
(\tilde g,\cA)_*:(C/V,\cA)^\sim\to(C/U,\cA)^\sim
$$
generalizing remark \ref{rem_continue-local}(i). We may then
consider the category
$$
(C/\bullet,\cA)^\sim
$$
whose objects are the pairs $(U,\cF)$ consisting of
objects $U\in\Ob(\cC)$ and $\cF\in\Ob(C/U,\cA)^\sim$.
A morphism $(U,\cF)\to(V,\cG)$ is the datum of a morphism
$g:U\to V$ of $\cC$ and a morphism $\cF\to(\tilde g,\cA)_*\cG$
in $(C/U,\cA)^\sim$. It follows easily from
\eqref{eq_restrict-sheaf} that the resulting functor
$$
(C/\bullet,\cA)^\sim\to\cC
\qquad
(U,\cF)\mapsto U
$$
is a fibration. Moreover, for every $U\in\Ob(\cC)$, every
covering family $(U_i\to U~|~i\in I)$ for $J_U$ generates
a sieve of $1$-descent for this fibration, which is even
of $2$-descent, if $\cA$ is complete and $\cC$ is small.
\end{remark}

\sset\subsubsection{}\label{subsec_Cat-star}
In this section we are mainly interested in presheaves and sheaves
with values in the category $\bCat$ of small categories, but the
following auxiliary construction shall also be useful. For every
category $\cB$ {\em whose fibred products are representable} we
consider the $2$-category
$$
\bCat^*(\cB)
$$
of {\em category objects} of $\cB$, whose objects are the data
$\cC^*:=(O,M,\ss,\st,\bone,\sc)$ such that $O,M\in\Ob(\cB)$,
and
\set\begin{equation}\label{eq_encode-Cats}
\xymatrix{ M \ar@<.5ex>[r]^-\ss \ar@<-.5ex>[r]_-\st & O}
\qquad
\bone:O\to M
\qquad
M\times_{(\st,\ss)}M\xrightarrow{\sc} M
\end{equation}
are morphisms in $\cB$ called respectively the {\em source},
{\em target}, {\em identity} and {\em composition laws} of $\cC^*$,
fulfilling the identities :
$$
\ss\circ\bone=\one_O=\st\circ\bone
$$
and making commute the diagrams of morphisms of $\cB$ :
$$
\xymatrix{ M\times_{(\st,\ss)}M\times_{(\st,\ss)}M
\ar[rr]^-{\one_M\times\sc} \ar[d]_{\sc\times\one_M} & &
M\times_{(\st,\ss)}M \ar[d]^\sc & 
M \ar[rr]^-{(\one_M,\bone\circ\st)} \ar[d]_{(\bone\circ\ss,\one_M)}
\ar[rrd]^-{\one_M} & &
M\times_{(\st,\ss)}M \ar[d]^\sc \\
M\times_{(\st,\ss)}M \ar[rr]^-\sc & & M &
M\times_{(\st,\ss)}M \ar[rr]^-\sc & & M.
}$$
Notice that the datum of $\sc$ includes the choice of a
representative for the fibre product $M\times_{(\st,\ss)}M$,
and we (implicitly) fix as well a universal cone
$M\leftarrow M\times_{(\st,\ss)}M\to M$ for this fibre product.

The $1$-cells $(F_1,F_2):\cC^*:=(O,M,\ss,\st,\bone,\sc)\to
\cC'^*:=(O',M',\ss',\st',\bone',\sc')$ in $\bCat^*(\cB)$ are
the pairs of morphisms $(F_1:O\to O',F_2:M\to M')$ of $\cB$
such that
$$
\ss'\circ F_2=F_1\circ\ss
\qquad
\st'\circ F_2=F_1\circ\st
\qquad
\bone'\circ F_1=F_2\circ\bone
\qquad
\sc'\circ(F_2\times_{(\st,\ss)}F_2)=F_2\circ\sc.
$$
The composition law of $1$-cells is given by the obvious rule :
$$
(F'_1,F'_2)\circ(F_1,F_2):=(F'_1\circ F_1,F'_2\circ F_2)
\quad
\text{for every pair of $1$-cells $\cC^*\xrightarrow{(F_1,F_2)}\cC'^*
\xrightarrow{(F'_1,F'_2)}\cC''^*$}.
$$
For every pair of category objects $\cC^*$ and $\cC'^*$ of $\cB$,
and every pair of $1$-cells $(F_1,F_2),(G_1,G_2):\cC^*\to\cC'^*$,
the $2$-cells $\beta:(F_1,F_2)\Rightarrow(G_1,G_2)$ are the morphisms
$\beta:O\to M'$ of $\cB$ such that
\set\begin{equation}\label{eq_spell-out-nat-transf}
\ss'\circ\beta=F_1
\qquad
\st'\circ\beta=G_1
\qquad
\sc'\circ(\beta\circ\ss,G_2)=\sc'\circ(G_1,\beta\circ\st)
\end{equation}
(where $(\beta\circ\ss,G_2)$ and $(G_1,\beta\circ\st)$ are
morphisms $M\to M'\times_{(\st',\ss')}M'$).

\begin{remark}\label{rem_Ob-and-Morph}
(i)\ \
We have a natural strict and strong $2$-equivalence
of $2$-categories :
\set\begin{equation}\label{eq_Cats-and-stars}
\sV\tdu\bCat\isom\bCat^*(\sV\tdu\Set)
\qquad
\text{for every universe $\sV$}.
\end{equation}
Namely, to each $\sV$-small category $\cC$ we assign the object
$[\cC]:=(\Ob(\cC),\rMorph(\cC),\ss,\st,\bone,\sc)$ where
$\ss,\st,\bone,\sc$ encode the source, target, identity and
composition laws for $\cC$ in the obvious way. Then any functor
$F:\cC\to\cC'$ induces an obvious $1$-cell $[F]:[\cC]\to[\cC']$,
and every natural transformation $\beta:F\Rightarrow F'$ induces
a $2$-cell $[\beta]:[F_1]\Rightarrow[F_2]$. The quasi-inverse
assigns to every datum $\cC^*:=(O,M,\ss,\st,\bone,\sc)$ the
category $[\cC^*]$ with $\Ob([\cC^*]):=O$ and
$\Hom_{[\cC^*]}(X,Y):=\ss^{-1}(X)\cap\st^{-1}(Y)$ for every
$X,Y\in O$, with the composition law induced by $\sc$ in
the obvious way, and with $\one_X:=\bone(X)$ for every $X\in O$.
Then every $1$-cell $(F_1,F_2):\cC^*\to\cC'^*$ induces a functor
$[F_1,F_2]:[\cC^*]\to[\cC'^*]$ and every $2$-cell
$\beta:(F_1,F_2)\Rightarrow(G_1,G_2)$ induces a natural
transformation $[\beta]:[F_1,F_2]\Rightarrow[G_1,G_2]$, in the
obvious fashion.

(ii)\ \
We have obvious {\em objects} and {\em morphisms} functors :
$$
\xymatrix{
\sV\tdu\bCat \ar@<.5ex>[rr]^-\Ob \ar@<-.5ex>[rr]_-\rMorph
& & \sV\tdu\Set
}\qquad
\text{for every universe $\sV$}
$$
which -- in terms of the $2$-equivalence
\eqref{eq_Cats-and-stars} -- translate as the functors
$\bCat^*(\sV\tdu\Set)\to\sV\tdu\Set$ that extract from each
datum $(O,M,\ss,\st,\bone,\sc)$ the set $O$ and respectively
the set $M$, and that assign to every morphism $(F_1,F_2)$ of
$\bCat^*(\sV\tdu\Set)$ the map $F_1$ and respectively the map $F_2$.

(iii)\ \
Notice that both $\Ob$ and $\rMorph$ are representable functors.
Indeed, let $\bone$ be the category with one object and one morphism,
and let $\mathbbm{2}$ be the category with exactly two objects $a$
and $b$, and a single morphism from $a$ to $b$ (and no morphisms
from $b$ to $a$). Then it is easily seen that $\bone$ (resp.
$\mathbbm{2}$) represents the presheaf $\Ob$ (resp. $\rMorph$)
on $\sV\tdu\bCat^o$ : the details shall be left to the reader.
Especially, $\Ob$ and $\rMorph$ commute with all the limits of
$\sV\tdu\bCat$ (example \ref{ex_lims-and-representables}(ii)).

Moreover, obviously a functor $F:\cA\to\cA'$ is an isomorphism
in $\sV\tdu\bCat$ if and only if both $\Ob(F):\Ob(\cA)\to\Ob(\cA')$
and $\rMorph(F):\rMorph(\cA)\to\rMorph(\cA')$ are bijections.

(iv)\ \
Taking into account (iii) and remark
\ref{rem_sheaves-with-values-in-A}(ii), it follows immediately
that a presheaf of $\sV$-small categories $F$ on $\cC$ is a
sheaf of categories on $C$ if and only if $\Ob(F):=\Ob\circ F$
and $\rMorph(F):=\rMorph\circ F$ are sheaves (of sets) on $C$.

(v)\ \
For every site $C:=(\cC,J)$, the categories
$(\cC,\sV\tdu\bCat)^\wedge$ and $(C,\sV\tdu\bCat)^\sim$ inherit
from $\sV\tdu\bCat$ a natural structure of $2$-category. Namely,
$(\cC,\sV\tdu\bCat)^\wedge$ can be described as the
sub-$2$-category of of $\sPsFun(\cC,\sV\tdu\bCat)$ whose objects
are the strict pseudo-functors, whose $1$-cells are all the strict
pseudo-natural transformations (where $\cC$ is regarded as usual,
as a $2$-category with trivial $2$-cells), and whose $2$-cells
are all the modifications $\beta\leadsto\beta'$, for every pair
of $1$-cells $\beta,\beta':F\Rightarrow F'$ of
$(\cC,\sV\tdu\bCat)^\wedge$. Similarly we describe
$(C,\sV\tdu\bCat)^\sim$ as the strong sub-$2$-category of
$(\cC,\sV\tdu\bCat)^\wedge$ whose objects are the sheaves of
$\sV$-small categories on the site $C$ (see definition
\ref{def_using-mods}(iii)). Then \eqref{eq_Cats-and-stars}
generalizes to strict and strong $2$-equivalences
\set\begin{equation}\label{eq_upgrade-Cats-and-stars}
(\cC,\sV\tdu\bCat)^\wedge\isom\bCat^*(\cC_\sV^\wedge)
\qquad
(C,\sV\tdu\bCat)^\sim\isom\bCat^*(C_\sV^\sim).
\end{equation}

(vi)\ \
Let $\cB$ and $\cB'$ be any two categories whose fibre products
are representable, and $u:\cB\to\cB'$ any functor {\em that
commutes with fibre products}. Then $u$ induces a strict
pseudo-functor
$$
\bCat^*(u):\bCat^*(\cB)\to\bCat^*(\cB')
\qquad
(O,M,\ss,\st,\bone,\sc)\mapsto(u(0),u(M),u(\ss),u(\bone),u(\sc)).
$$
Indeed, we can regard $u(\sc)$ as a morphism
$u(M)\times_{(u(\st),u(\ss))}u(M)\to u(M)$ (since $u$ commutes with
fibre products), and then it is clear that the rule defining
$\bCat^*(u)$ takes objects of $\bCat^*(\cB)$ to objects
of $\bCat^*(\cB')$. To every $1$-cell $(F_1,F_2)$ and every
$2$-cell $\beta$, the pseudo-functor $\bCat^*(u)$ assigns
likewise respectively the $1$-cell $(u(F_1),u(F_2))$ and the
$2$-cell $u(\beta)$ of $\bCat^*(\cB')$.

If $v:\cB\to\cB'$ is another functor commuting with fibre products,
and $\alpha:u\Rightarrow v$ is any natural transformation, we get
a strict pseudo-natural transformation
$$
\bCat^*(\alpha):\bCat^*(u)\Rightarrow\bCat^*(v)
\qquad
(O,M,\ss,\st,\bone,\sc)\mapsto(\alpha_O,\alpha_M).
$$
Indeed, notice that $u(M\times_{(\st,\ss)}M)$ represents
$u(M)\times_{(u(\st),u(\ss))}u(M)$ since $u$ commutes with fibre
products, and likewise for $v(M\times_{(\st,\ss)}M)$, and under
these identifications, the morphism $\alpha_{M\times_{(\st,\ss)}M}$
corresponds to $\alpha_M\times_{(u(\st),u(\ss))}\alpha_M$, whence
it is easily seen that $(\alpha_O,\alpha_M)$ is a $1$-cell of
$\bCat^*(\cB')$, and the strict pseudo-functoriality of
$\bCat^*(\alpha)$ is then immediate from the definition.

(vii)\ \
Clearly the rules $u\mapsto\bCat^*(u)$ and
$\alpha\mapsto\bCat^*(\alpha)$ define a strict pseudo-functor
from the sub-$2$-category of $\bCat$ whose objects are all
the small categories whose fibre products are representable,
whose $1$-cells are the functor commuting with fibre products
and whose $2$-cells are the natural transformations, to the
$2$-category whose objects are the small $2$-categories, whose
$1$-cells are the strict pseudo-functors, and whose $2$-cells
are the strict pseudo-natural transformations.
\end{remark}

\sset\subsubsection{}\label{subsec_sheafify-categories}
Let now $C:=(\cC,J)$ be a $\sU$-site; then we deduce that the
strongly faithful inclusion pseudo-functor
$(C,\sV\tdu\bCat)^\sim\to(\cC,\sV\tdu\bCat)^\wedge$ admits a
strong left $2$-adjoint :
\set\begin{equation}\label{eq_sheafify-on-Cat}
(\cC,\sV\tdu\bCat)^\wedge\to(C,\sV\tdu\bCat)^\sim
\qquad
F\mapsto F^a
\end{equation}
for every universe $\sV$ with $\sU\subset\sV$. Namely, by
remark \ref{rem_Ob-and-Morph}(vi,vii), the exact left adjoint
$(-)^a:\cC^\wedge_\sV\to C^\sim_\sV$ provided by remark
\ref{rem_U-site}(ii) induces first a strong left $2$-adjoint
$\bCat^*((-)^a):\bCat^*(\cC_\sV^\wedge)\to\bCat^*(C_\sV^\sim)$ to the
inclusion of $2$-categories $\bCat^*(C^\sim_\sV)\to\bCat^*(C^\wedge_\sV)$.
Then, combining with the strong and strict $2$-equivalences
\eqref{eq_upgrade-Cats-and-stars} we get the sought strong left
$2$-adjoint. It is also easily seen that \eqref{eq_sheafify-on-Cat}
is an exact functor on the underlying categories.

\sset\subsubsection{}\label{subsec_pull-back-for-Cats}
Let $C:=(\cC,J)$ and $C':=(\cC',J')$ be two $\sU$-sites,
$u:\cC\to\cC'$ a functor, $\sU'$ a universe with $\sU\subset\sU'$
and such that $\cC$ and $\cC'$ are $\sU'$-small and the functor
$u^a_{\sU'!}:\cC^\wedge_{\sU'}\to C'^\sim_{\sU'}$ commutes with fibre
products. Arguing as in \eqref{subsec_sheafify-categories}, and
in light of remark \ref{rem_change-universe}(iii), we deduce that
for every universe $\sV$ with $\sU\subset\sV$, the pseudo-functor
$u^a_{\sV!}$ induces a pseudo-functor
$$
(u,\sV\tdu\bCat)^a_!:(\cC,\sV\tdu\bCat)^\wedge\to(C',\sV\tdu\bCat)^\sim.
$$
In case $u$ is continuous for the topologies $J$ and $J'$,
corollary \ref{cor_two-U-sites}(i) and remark
\ref{rem_Ob-and-Morph}(vi,vii) easily imply that $(u,\sV\tdu\bCat)^a_!$
restricts to a strong and strict left $2$-adjoint
$$
(\tilde u,\sV\tdu\bCat)^*:(C,\sV\tdu\bCat)^\sim\to(C',\sV\tdu\bCat)^\sim
$$
for the strict pseudo-functor $(\tilde u,\sV\tdu\bCat)_*$ of remark
\ref{rem_sheaves-with-values-in-A}(iv).

\sset\subsubsection{}
Let $C:=(\cC,J)$ and $C':=(\cC',J')$ be two $\sU$-sites,
$u:\cC\to\cC'$ a cocontinuous functor, and $\sV$ a universe
with $\sU\subset\sV$. Arguing as in
\eqref{subsec_sheafify-categories} we see that the adjoint
pair of left exact functors $(\breve u{}^*_\sV,\breve u_{\sV*})$
provided by corollary \ref{cor_two-U-sites}(ii) induces a
strict and strong $2$-adjoint pair of strict pseudo-functors
$$
\xymatrix@C+30pt{ (C',\sV\tdu\bCat)^\sim
\ar@<.5ex>[r]^-{(\breve u,\sV\tdu\bCat)^*} & (C,\sV\tdu\bCat)^\sim.
\ar@<.5ex>[l]^-{(\breve u,\sV\tdu\bCat)_*}
}$$

\begin{remark}\label{rem_alternative-inverse-image}
(i)\ \
In the situation of \eqref{subsec_pull-back-for-Cats}, suppose
that $u$ is a continuous functor and $\tilde u_{\sV*}$ is an
equivalence, so that $\tilde u{}_\sV^*$ is its quasi-inverse. From
the explicit description of the functors $(\tilde u,\sV\tdu\bCat)_*$
and $(\tilde u,\sV\tdu\bCat)^*$, we then deduce straightforwardly
that the latter are mutually quasi-inverse equivalences and strong
$2$-equivalences as well.

(ii)\ \
Combining with theorem \ref{th_canon-topos}(iii), we deduce
especially that for every $\sU$-site $C:=(\cC,J)$, the Yoneda
embedding $h^a_\cC:\cC\to C_\sU^\sim$ induces a strong $2$-equivalence
of $2$-categories :
$$
(\tilde h{}^a_\cC,\sV\tdu\bCat)_*:
(\Can(C_\sU^\sim),\sV\tdu\bCat)^\sim\isom(C,\sV\tdu\bCat)^\sim.
$$

(iii)\ \
Let $C:=(\cC,J)$ and $C':=(\cC',J')$ be two sites such that $\cC$
is small and $\cC'$ has small $\Hom$-sets, and let $u:\cC\to\cC'$
be any continuous functor for the topologies $J$ and $J'$. Recall
that for every universe $\sV$ with $\sU\subset\sV$, the functor
$(u,\sV\tdu\bCat)^\wedge$ admits a left adjoint
$$
(u,\sV\tdu\bCat)_!:(\cC,\sV\tdu\bCat)^\wedge\to(\cC',\sV\tdu\bCat)^\wedge
$$
computed by left Kan extensions : see theorem \ref{th_Kan-ext}.
If $C'$ is a $\sU$-site, we can compose this functor with the
functor \eqref{eq_sheafify-on-Cat}, to get a functor
$(u,\sV\tdu\bCat)^a_!:(\cC,\sV\tdu\bCat)^\wedge\to(C',\sV\tdu\bCat)^\sim$,
whose restriction to $(\cC,\sV\tdu\bCat)^\sim$ is a left adjoint for
$(\tilde u,\sV\tdu\bCat)_*$. If $C$ is no longer a small site, but
only a $\sU$-site, we can pick as usual a small topologically
generating full subcategory $\cG\subset\cC$, endowed with
the topology $J_\cG$ induced by $J$ via the inclusion functor
$j:\cG\to\cC$, so that
$(\tilde\jmath,\sV\tdu\bCat)_*:(C,\sV\tdu\bCat)^\sim\to
((\cG,J_\cG),\sV\tdu\bCat)^\sim$ is an equivalence, by (i); then
$(\tilde u,\sV\tdu\bCat)_*$ admits the left adjoint
$(\widetilde{u\circ\jmath},\sV\tdu\bCat)^*\circ
(\tilde\jmath,\sV\tdu\bCat)_*$.
This yields another construction for the functor underlying the
strict pseudo-functor $(\tilde u,\sV\tdu\bCat)^*$ of
\eqref{subsec_pull-back-for-Cats}, which is available under
weaker assumptions, but is not obviously pseudo-functorial for
the $2$-category structures of $(C,\sV\tdu\bCat)^\sim$ and
$(C',\sV\tdu\bCat)^\sim$.

(iv)\ \
On the other hand, in the situation of (iii), recall that
the strict pseudo-functor
$$
\sPsFun(u^o,\sV\tdu\bCat):\sPsFun(\cC'^o,\sV\tdu\bCat)\to
\sPsFun(\cC^o,\sV\tdu\bCat)
$$
admits a left $2$-adjoint that can be computed by the strong
left $2$-Kan extension $2\tdu\!\!\int^{u^o}$ along $u^o$ (remark
\ref{rem_strong-2-Kan-ext}(ii) and theorem \ref{th_bCat-cpt}).
The latter restricts to a strict pseudo-functor
$$
(u,\sV\tdu\bCat)_!:
(\cC,\sV\tdu\bCat)^\wedge\to(\cC',\sV\tdu\bCat)^\wedge.
$$
Moreover, the $2$-adjunction for the pair
$(2\tdu\!\!\int^{u^o},\sPsFun(u^o,\sV\tdu\bCat))$ restricts
to a natural functor
$$
(-)^\dagger_{F,G}:\Hom_{(\cC,\sV\tdu\bCat)^\wedge}(G\circ u^o,F)\to
\Hom_{(\cC,\sV\tdu\bCat)^\wedge}\Bigl(G,2\tdu\!\!\!\int^{u^o} F\Bigr)
$$
for every presheaf of categories $F$ on $\cC$ and $G$ on
$\cC'$ (remark \ref{rem_cont-strong-2-Kan-ext}). But in
this generality, though such functors are fully faithful,
it is not clear whether they are equivalences of categories.

(v)\ \
However, {\em suppose now additionally that the category
$X/u\cC$ is cofiltered for every $X\in\Ob(\cC')$} (this
condition is fulfilled, notably, in case {\em $C$ and $C'$ are
two lex-sites, and $u$ is a morphism of lex-sites $C'\to C$},
by virtue of example \ref{ex_cofiltered-comma}(i)).
Recall that for every presheaf of $\sV$-small categories $F$ on
$\cC$ and every $X\in\Ob(\cC')$, the value $2\tdu\!\!\int^{u^o} F(X)$
represents the $2$-colimit of the strict pseudo-functor
$F\circ\st^o_X:(X/u\cC)^o\to\sV\tdu\bCat$; but then such
$2$-colimit is represented by the colimit of the same functor,
and every universal cocone for such colimit is also a universal
pseudo-cocone (example \ref{ex_filter-2-colim-in-Cat}(iv)).
This means that the restriction to $(\cC,\sV\tdu\bCat)^\wedge$
of the strong left $2$-Kan extension $2\tdu\!\!\int^{u^o}$ agrees
with the left Kan extension $\int^{u^o}$ of theorem \ref{th_Kan-ext};
moreover, a direct inspection shows that the functor $(-)^\dagger_{F,G}$
agrees on objects with the standard adjunction for the pair
$(\int^{u^o},(u,\sV\tdu\bCat)^\wedge)$, provided by the proof of
theorem \ref{th_Kan-ext}. {\em I.e.}, under these assumptions,
$(u,\sV\tdu\bCat)_!$ is indeed a strong and strict left
$2$-adjoint for $(u,\sV\tdu\bCat)^\wedge$, and we get an
essentially commutative diagram of $2$-categories :
$$
\xymatrix@C+20pt{
(\cC,\sV\tdu\bCat)^\wedge \ar[r]^-{(u,\sV\tdu\bCat)_!} \ar[d] &
(\cC',\sV\tdu\bCat)^\wedge \ar[d] \\
\sPsFun(\cC^o,\sV\tdu\bCat) \ar[r]^-{2\tdu\!\!\int^{u^o}} &
\sPsFun(\cC'^o,\sV\tdu\bCat)
}$$
whose vertical arrows are the inclusion strict pseudo-functors
(with a different choice of left $2$-Kan extension, this diagram
is then still at least pseudo-commutative). Lastly, under these
same assumptions, we obtain as well a pseudo-commutative diagram :
$$
\xymatrix@C+30pt{
(\cC,\sV\tdu\bCat)^\wedge \ar[r]^-{(u,\sV\tdu\bCat)_!}
\ar[d]_{\cFib_\cC} & (\cC',\sV\tdu\bCat)^\wedge \ar[d]^{\cFib_{\cC'}} \\
\sV\tdu\Fib(\cC) \ar[r]^-{\sV\tdu\Fib(u)_!} & \sV\tdu\Fib(\cC').
}$$
Arguing as in (iii), we get, again under the current assumptions,
an alternative construction of the strict pseudo-functor
$(\tilde u,\sV\tdu\bCat)^*$, as the restriction to
$(\cC,\sV\tdu\bCat)^\sim$ of the composition of the foregoing
pseudo-functor $(u,\sV\tdu\bCat)_!$ with the pseudo-functor
\eqref{eq_sheafify-on-Cat}.
\end{remark}

\begin{lemma}\label{lem_stacks-and-sheaves-of-cats}
Let $C:=(\cC,J)$ be any site, and consider a presheaf of
categories on $\cC$
$$
\cA_\bullet:\cC^o\to\bCat
\qquad
X\mapsto\cA_X
\qquad
(f:Y\to X)\mapsto(\cA_f:\cA_X\to\cA_Y).
$$
The following holds :
\begin{enumerate}
\item
If $\cA_\bullet$ is a sheaf on $C$, the associated prestack
$\cFib(\cA_\bullet)$ on $C$ is $1$-separated.
\item
Suppose that $C$ is a $\sU$-site, and that $\cFib(\cA_\bullet)$ is
a stack on $C$, and let $j:\cA_\bullet\to\cA^a_\bullet$ be the unit
of adjunction (notation of \eqref{subsec_sheafify-categories}).
Then the induced morphism of prestacks
$\cFib(j):\cFib(\cA_\bullet)\to\cFib(\cA^a_\bullet)$ is an
equivalence of categories.
\end{enumerate}
\end{lemma}
\begin{proof}(i): By lemma \ref{lem_crit-0-1_separation}, it
suffices to show that for every $X\in\Ob(\cC)$ and every pair
of cartesian sections $\sigma,\sigma':\cC/X\to\cFib(\cA_\bullet)$,
the presheaf $\cCart(\sigma,\sigma')$ is a sheaf on the site $C/X$.
Set $A:=\sigma(\one_X)$ and $A':=\sigma'(\one_X)$; it is easily
seen that $\cCart(\sigma,\sigma')$ is isomorphic to the presheaf
$$
H_{AA'}:(\cC/X)^o\to\Set
\qquad
(f:Y\to X)\mapsto\Hom_{\cA_Y}(\cA_f(A),\cA_f(A'))
$$
that assigns to every morphism
$(h/X):(Y\xrightarrow{f}X)\to(Y'\xrightarrow{f'}X)$ of $\cC/X$
the induced map
$$
\Hom_{\cA_Y'}(\cA_{f'}(A),\cA_{f'}(A'))\to\Hom_{\cA_Y}(\cA_f(A),\cA_f(A'))
\qquad
g\mapsto\cA_h(g).
$$
Now, let $\bone_{\cC/X}$ be a final object of $(\cC/X)^\wedge$ (see
example \ref{ex_equalizers}(v)); it is easily seen that
$\bone_{\cC/X}$ is a sheaf on $C/X$. Recall that the source
functor $\ss_X:\cC/X\to\cC$ is continuous for the topologies of $C$
and $C/X$. Denote also by $(O_\cA,M_\cA,\ss_\cA,\st_\cA,\bone_\cA)$
the object of $(\cC,\bCat^*)^\wedge$ corresponding to $\cA_\bullet$;
by remark \ref{rem_Ob-and-Morph}(iv), both $O_\cA$ and $M_\cA$
are sheaves on $C$, hence we get sheaves $\tilde\ss_{X*}(O_\cA)$
and $\tilde\ss_{X*}(M_\cA)$ on $C/X$, as well as source and target
morphisms of sheaves
$\tilde\ss_{X*}(\ss_\cA),\tilde\ss_{X*}(\st_\cA):\tilde\ss_{X*}(M_\cA)
\to\tilde\ss_{X*}(O_\cA)$. There follows a cartesian diagram of
presheaves :
$$
{\diagram H_{AA'} \ar[r] \ar[d] & \tilde\ss_{X*}(M_\cA)
\ar[d]^{(\tilde\ss_{X*}(\ss_\cA),\tilde\ss_{X*}(\st_\cA))} \\
\bone_{\cC/X} \ar[r]^-\tau & \tilde\ss_{X*}(O_\cA\times O_\cA)
\enddiagram}
\qquad
\text{ on $\cC/X$}$$
where $\tau_f:\{\bone\}=\bone_{\cC/X}(X)\to\Ob(\cA_Y\times\cA_Y)$
is the map such that $\bone\mapsto(\cA_f(A),\cA_f(A'))$, for every
morphism $f:Y\to X$ of $\cC$. But then $H_{AA'}$ is a sheaf on $C/X$,
as required.

(ii): By assumption, $\cFib(\cA_\bullet)$ is $2$-separated, and
$\cFib(\cA^a_\bullet)$ is $1$-separated, by (i); in view of
proposition \ref{prop_faith-separ-cover}, it then suffices
to check that $\cFib(j)$ is $i$-covering for $i=0,1,2$. The
latter follows easily from remark \ref{rem_realign}, taking
into account the explicit description of $\cA^a_\bullet$ from
\eqref{subsec_sheafify-categories}.
\end{proof}

\sset\subsubsection{}\label{subsec_ripristine}
\label{subsec_u-with-ripristine}
Fix a universe $\sV$ with $\sU\!\subset\!\sV$, and a $\sU$-site
$C:=(\cC,J_\cC)$. Consider the pseudo-functor
$$
\cFib^a_C:(C,\sV\tdu\bCat)^\sim\to\sV\tdu\Stack(C)
\qquad
\cA_\bullet\mapsto\cFib(\cA_\bullet)^a
$$
composition of the restriction of the pseudo-functor $\cFib_\cC$
of theorem \ref{th_fundamental-fibrations} with $(-)^a$ of corollary
\ref{cor_stackificate-U-sites} (where $(C,\sV\tdu\bCat)^\sim$ is a
$2$-category as in remark \ref{rem_Ob-and-Morph}(v)).
Let also $C':=(\cC',J_{\cC'})$ be another $\sU$-site
and $u:C'\to C$ a weak morphism of sites (definition
\ref{def_weak-morph-of-sites}); we wish to attach an
orientation to the induced diagram :
\set\begin{equation}\label{eq_oriented-by-Delta-u}
{\diagram (C',\sV\tdu\bCat)^\sim \ar[rr]^-{(\tilde u,\sV\tdu\bCat)_*}
\ar[d]_{\cFib^a_{C'}} & & (C,\sV\tdu\bCat)^\sim \ar[d]^{\cFib^a_C} \\
\sV\tdu\Stack(C') \ar[rr]^-{\sV\tdu\St(u)_*} & & \sV\tdu\Stack(C)
\enddiagram}\end{equation}
{\em i.e.} a pseudo-natural transformation
$$
\Delta^u:\cFib^a_C\circ(\tilde u,\sV\tdu\bCat)_*\Rightarrow
\sV\tdu\St(u)_*\circ\cFib^a_{C'}.
$$
To this aim, notice that remark
\ref{rem_distinguished-cleavage}(ii) yields an oriented
diagram of $2$-categories :
\set\begin{equation}\label{eq_with-perps}
{\spreaddiagramcolumns{+30pt}\diagram
(C',\sV\tdu\bCat)^\sim \ar[r]^-{(\tilde u,\sV\tdu\bCat)_*}
\ar[d]_{\cFib_{\cC'}} \drtwocell\omit{_\ \ \ \perp^u} &
(C,\sV\tdu\bCat)^\sim \ar[d]^{\cFib_\cC} \\
\sV\tdu\Fib(\cC') \ar[r]_-{\sV\tdu\Fib(u)^*} & \sV\tdu\Fib(\cC)
\enddiagram}\end{equation}
whose orientation $\perp^u$ is a strict pseudo-natural
isomorphism of strict pseudo-functors. On the other hand,
the square of links $\cD_u$ of
\eqref{subsec_weak-morph-of-sites-and-cocont}, oriented by
the identity pseudo-natural transformation, yields a base
change square :
\set\begin{equation}\label{eq_oriented-by-Delta_1-u}
{\spreaddiagramcolumns{+30pt}\diagram
\sV\tdu\Fib(\cC') \ar[r]^-{\sV\tdu\Fib(u)^*}
\ar[d]_{(-)_{C'}^a} \drtwocell\omit{_\ \ \ \ \ \ \ \Upsilon(\cD_u)} &
\sV\tdu\Fib(\cC) \ar[d]^{(-)^a_C} \\
\sV\tdu\Stack(C') \ar[r]_-{\sV\tdu\St(u)_*} & \sV\tdu\Stack(C)
\enddiagram}
\end{equation}
and then we take $\Delta^u:=\perp^u\boxvert\Upsilon(\cD_u)$.
In case we wish to emphasize the dependence on $\sV$, we
shall also write $\sV\tdu\Delta^u$ for this orientation of
\eqref{eq_oriented-by-Delta-u}.

\begin{remark}\label{rem_get-link-if-u-strong}
In the situation of \eqref{subsec_u-with-ripristine}, suppose
that $u$ is a morphism of sites; in this case, the
pseudo-functors $(\tilde u,\sV\tdu\bCat)_*$ and
$\sV\tdu\St(u)_*$ admit left $2$-adjoints
$(\tilde u,\sV\tdu\bCat)^*$ and $\sV\tdu\St(u)^*$, and after
choosing units and counits for these $2$-adjoint pairs, and
for the $2$-adjoint pair $(\sV\tdu\Fib(u)_!,\sV\tdu\Fib(u)^*)$
we then get well defined links
$$
\begin{aligned}
(\tilde u,\sV\tdu\bCat):=
&\,((\tilde u,\sV\tdu\bCat)^*,(\tilde u,\sV\tdu\bCat)_*,
\eta^{(\tilde u,\sV\tdu\bCat)},\eps^{(\tilde u,\sV\tdu\bCat)})\!:\!
(C',\sV\tdu\bCat)^\sim\!\!\to\!(C,\sV\tdu\bCat)^\sim \\
\sV\tdu\St(u):=&\,(\sV\tdu\St(u)^*,\sV\tdu\St(u)_*,
\eta^{\sV\tdu\St(u)},\eps^{\sV\tdu\St(u)}):
\sV\tdu\Stack(C')\to\sV\tdu\Stack(C) \\
\sV\tdu\Fib(u):=&(\sV\tdu\Fib(u)_!,\sV\tdu\Fib(u)^*,
\eta^{\sV\tdu\Fib(u)},\eps^{\sV\tdu\Fib(u)}):
\sV\tdu\Fib(\cC')\to\sV\tdu\Fib(\cC)
\end{aligned}
$$
of the $2$-category $\sV\tdu\overline{2\tdu\bCat}$,
and the data
$$
((-)^a_{C'},(-)^a_C,\Upsilon(\cD_u))
\qquad\text{and}\qquad
(\cFib_{\cC'},\cFib_\cC,\perp^u)
$$
can be regarded as $1$-cells $\sV\tdu\Fib(u)\to\sV\tdu\St(u)$
and $(\tilde u,\sV\tdu\bCat)\to\sV\tdu\Fib(u)$ in the $2$-category
$\wLink(\sV\tdu\overline{2\tdu\bCat})$ (see
\eqref{subsec_variant-of-Upsilon}) whose composition is the $1$-cell
$$
(\cFib_{C'}^a,\cFib^a_C,\Delta^u):
(\tilde u,\sV\tdu\bCat)\to\sV\tdu\St(u).
$$
\end{remark}

\sset\subsubsection{}\label{subsec_Fib-Stack-Cat}
Let $\sV$ be a universe with $\sU\in\sV$; we shall denote by
$$
(\sU,\sV)\tdu\wSite
$$
the sub-$2$-category of $(\sU,\sV)\tdu\Site$ (see definition
\ref{def_morph-of-sites}(ii)) whose objects are all the
$\sV$-small $\sU$-sites, whose $1$-cells are the weak morphisms
of sites, and whose $2$-cells $g\Rightarrow g'$ are the natural
transformations of functors $g'\Rightarrow g$.
As in remark \ref{rem_morph-of-sites}(i), we shall usually
drop the mention of $\sU$ and $\sV$, and write simply $\wSite$
for this $2$-category. We have a strict pseudo-functor
$$
(-,\sV\tdu\bCat)^\sim:\wSite\to\sV\tdu\overline{2\tdu\bCat}
\qquad
C\mapsto(C,\sV\tdu\bCat)^\sim
$$
(notation of remark \ref{rem_reduced-2-cats}) that assigns to
every morphism of $\sU$-sites $u$ the strict pseudo-functor
$(\tilde u,\sV\tdu\bCat)_*$, and to every natural
transformation $\beta:u\Rightarrow v$ the $2$-cell
$(\tilde\beta,\sV\tdu\bCat)_*:(\tilde v,\sV\tdu\bCat)_*
\Rightarrow(\tilde u,\sV\tdu\bCat)_*$. Likewise, we get
a pseudo-functor
$$
\sV\tdu\Stack:\wSite\to\sV\tdu\overline{2\tdu\bCat}
\qquad
C\mapsto\sV\tdu\Stack(C)
\qquad
(u:C'\to C)\mapsto\sV\tdu\St(u)_*.
$$
To define the coherence constraints
$(\delta^{\sV\tdu\Stack},\gamma^{\sV\tdu\Stack})$
of $\sV\tdu\Stack$, and the action of $\sV\tdu\Stack$ on natural
transformations $\beta:u\Rightarrow v$ between weak morphisms of
sites $u,v:C'\to C$, notice that
$\sV\tdu\Stack(\cC,J)\subset\sV\tdu\Fib(\cC)$ for every
$\sU$-site $(\cC,J)$, and $\sV\tdu\St(u)_*$ is the restriction
of $\sV\tdu\Fib(u)^*$, for every morphism of sites $u$; we may
then simply define $(\delta^{\sV\tdu\Stack},\gamma^{\sV\tdu\Stack})$
as the restriction of the coherence constraint
$(\delta^{\sV\tdu\Fib},\gamma^{\sV\tdu\Fib})$ of the pseudo-functor
$\sV\tdu\Fib$ of \eqref{subsec_complete-Fib}, and for any $\beta$
as above, let
$\sV\tdu\St(\beta)_*:\sV\tdu\St(v)\Rightarrow\sV\tdu\St(u)_*$ be
the restriction of $\sV\tdu\Fib(\beta)^*$.

With this notation, it would be tempting to state that the
rule : $C\mapsto\cFib^a_C$ of \eqref{subsec_u-with-ripristine}
defines a pseudo-natural transformation
$\cFib^a_\bullet:(-,\sV\tdu\bCat)^\sim\Rightarrow\sV\tdu\Stack$
whose coherence constraint is given by the orientations
$\Delta^u$. However, notice first that $\Delta^u$ points
in the direction opposite to the one which is required for
such a coherence constraint. This can be fixed, by stating
instead that $\cFib^a_\bullet$ should be a pseudo-natural
transformation
${}^o(-,\sV\tdu\bCat)^\sim\Rightarrow{}^o\sV\tdu\Stack$. But
since {\em the orientations $\Delta^u$ are not, in general,
pseudo-natural equivalences}, we have rather :

\begin{proposition}\label{prop_Brexit}
There exists a lax-natural transformation
$$
{}^o\cFib^a_\bullet:
{}^o(-,\sV\tdu\bCat)^\sim\Rightarrow{}^o\sV\tdu\Stack
\qquad
C\mapsto{}^o\cFib^a_C
\qquad
\text{for every $\sU$-site $C$}
$$
whose coherence constraint attaches to every weak morphism of\/
$\sU$-sites $u:C'\to C$ the square ${}^o\eqref{eq_oriented-by-Delta-u}$,
with its orientation ${}^o\Delta^u$.
\end{proposition}
\begin{proof} Consider the forgetful strict pseudo-functor
$$
\Phi:\wSite\to{}^o\sV\tdu\bCat^o
\qquad
(\cC,J)\mapsto\cC
$$
that assigns to every weak morphism of sites $u:(\cC',J')\to(C,J)$
the functor $u:\cC\to\cC'$, and to every $2$-cell $\beta:u\Rightarrow v$
of $\wSite$ the natural transformation $\beta:v\Rightarrow u$.
From corollary \ref{cor_murder-in-nice} we deduce a pseudo-natural
transformation of pseudo-functors :
$$
{}^o\cFib_\bullet:
{}^o(-,\sV\tdu\bCat)^\sim\Rightarrow{}^o\sV\tdu\Fib\circ{}^o\Phi
\qquad
(\cC,J)\mapsto{}^o\cFib_\cC
\qquad
\text{for every $\sU$-site $(\cC,J)$}
$$
whose coherence constraint assigns to every morphism of
sites $u$ the oriented square diagram \eqref{eq_with-perps}.
It then suffices to exhibit a lax-natural transformation
$$
{}^o(-)^a_\bullet:
{}^o\sV\tdu\Fib\circ{}^o\Phi\Rightarrow{}^o\sV\tdu\Stack
\qquad
C\mapsto{}^o(-)^a_C
\qquad
\text{for every $\sU$-site $C$}
$$
whose coherence constraint assigns to every weak morphism
of\/ $\sU$-sites $u:C'\to C$ the oriented square $\Upsilon(\cD_u)$.
Indeed, we will then define
${}^o\cFib^a_\bullet:={}^o(-)^a_\bullet\odot{}^o\cFib_\bullet$.
Now, notice that the pseudo-functors $\sV\tdu\Stack$ and
$\sV\tdu\Fib$ factor through well-defined pseudo-functors
$$
\sV\tdu\LStack:\wSite\to\sLink(\sV\tdu\overline{2\tdu\bCat})
\qquad
\sV\tdu\LFib:{}^o\sV\tdu\bCat^o\to\sLink(\sV\tdu\overline{2\tdu\bCat}).
$$
Namely, $\sV\tdu\LStack$ assigns to every site $C$ the $2$-category
$\sV\tdu\Stack(C)$ and to every weak morphism of sites $u:C'\to C$
the link $\sV\tdu\St(u)$ defined as in remark
\ref{rem_get-link-if-u-strong}; likewise, $\sV\tdu\LFib$ assigns
to every category $\cC$ the $2$-category $\sV\tdu\Fib(\cC)$, and
to every functor $u:\cC\to\cC'$ the link $\sV\tdu\Fib(u)$.
Next, let $i_C:\sV\tdu\Stack(C)\to\sV\tdu\Fib(\cC)$ be the inclusion
strict pseudo-functor; notice that the rule : $C\mapsto i_C$
for every site $C$ yields a pseudo-natural transformation
$$
i_\bullet:\sV\tdu\LStack\Rightarrow\sV\tdu\LFib\circ\Phi
$$
whose coherence constraints are given by the oriented
diagrams $\cD_u$ (notice that $i_\bullet$ is not strict,
even though $\cD_u$ is oriented by the identity pseudo-natural
transformation). According to remark \ref{rem_pseudo-natural}(ii),
the pseudo-natural transformation $i_\bullet$ corresponds to a
pseudo-functor $\tilde\imath_\bullet:
\wSite\to\cM:=2\tdu\sMorph(\sLink(\sV\tdu\overline{2\tdu\bCat}))$,
and we consider the composition :
$$
\Psi:\wSite\xrightarrow{\tilde\imath_\bullet}\cM
\xrightarrow{\Upsilon^o_{\sV\tdu\overline{2\tdu\bCat}}}
(2\tdu\sMorph(\sV\tdu\overline{2\tdu\bCat}^o))^o\isom
{}^o2\tdu\sMorph({}^o\sV\tdu\overline{2\tdu\bCat})
$$
with the pseudo-functor $\Upsilon^o_{\sV\tdu\overline{2\tdu\bCat}}$
of theorem \ref{th_Upsilon-is-pseudo-fctr} and the strict
isomorphism of $2$-categories given by example
\ref{ex_change-orientation}(i). Let
$\ss,\st:2\tdu\sMorph({}^o\sV\tdu\overline{2\tdu\bCat})\to
{}^o\sV\tdu\overline{2\tdu\bCat}$ be the source and target
strict pseudo-functors; by definition, ${}^o\Psi$ is a
lax-natural transformation
$\ss\circ{}^o\Psi\Rightarrow\st\circ{}^o\Psi\Rightarrow$.
Lastly, a direct inspection shows that
$\ss\circ{}^o\Psi=\LStack$ and
$\st\circ{}^o\Psi\Rightarrow=\LFib\circ\Phi$, so we may let
${}^o(-)^a_\bullet:={}^o\Psi$.
\end{proof}

\sset\subsubsection{}\label{subsec_rotate-the-cube}
We consider now an oriented square :
$$
\xymatrix{ D':=(\cD',J_{\cD'}) \ar[r]^-{v'} \ar[d]_{u'}
\drtwocell\omit{_\ \beta} & C':=(\cC',J_{\cC'}) \ar[d]^u \\
D:=(\cD,J_\cD) \ar[r]^-v & C:=(\cC,J_\cC)
}$$
where $v$ and $v'$ are morphisms of $\sU$-sites, and $u$
and $u'$ are weak morphisms of $\sU$-sites. In case also
$u$ and $u'$ are morphisms of sites, we get an induced
oriented diagram of topoi :
\set\begin{equation}\label{eq_links-from-topoi}
{\spreaddiagramcolumns{+10pt}\diagram
\tilde D{}' \ar[r]^-{\tilde v'} \ar[d]_{\tilde u'}
\drtwocell\omit{_\ \ \ \tilde\beta_*} & \tilde C' \ar[d]^{\tilde u} \\
\tilde D \ar[r]^-{\tilde v} & \tilde C.
\enddiagram}
\end{equation}
Even when $u$ and $u'$ are only weak morphisms of sites,
both $\tilde u_*$ and $\tilde u'_*$ admit left adjoint
functors, hence \eqref{eq_links-from-topoi} may still be
regarded as an oriented square of links in the $2$-category
$\sV\tdu\bCat$.  On the other hand, we deduce a diagram of
$2$-categories :
$$
\xymatrix@C+20pt{ \sV\tdu\Stack(D') \ar[rrr]^-{\St(v')_*}
\ar[ddd]_{\St(u')_*} & &
\dltwocell\omit{_\ \ \ \ \ \ \ \ \ \ \ \ \Delta^{v'}} &
\sV\tdu\Stack(C') \ar[ddd]^{\St(u)_*} \\
& \dltwocell\omit{_\ \Delta^{u'}\ \ \ \ \ }
(D',\sV\tdu\bCat)^\sim \ar[r]^-{(\tilde v',\sV\tdu\bCat)_*} \ar[lu]
\ar[d]_{(\tilde u',\sV\tdu\bCat)_*}
\drtwocell\omit{_\ \ \ \ \ \ \ \ \ \ \ \ (\tilde\beta,\sV\tdu\bCat)_*} &
(C',\bCat)^\sim \ar[d]^{(\tilde u,\sV\tdu\bCat)_*} \ar[ru]
\drtwocell\omit{^\ \ \ \ \ \ \ \ \ \ \ \ \Delta^u} \\
& (D,\sV\tdu\bCat)^\sim \ar[r]_-{(\tilde v,\sV\tdu\bCat)_*} \ar[ld]
\drtwocell\omit{_\ \Delta^v\ \ \ \ \ \ \ \ \ \ \ \ } &
(C,\sV\tdu\bCat)^\sim \ar[rd] & \\
\sV\tdu\Stack(D) \ar[rrr]_-{\St(v)_*} & & & \sV\tdu\Stack(C)
}$$
whose diagonal arrows are the pseudo-functors $\cFib^a$.
We complete it by adding the orientation
$$
\St(\beta)^\gamma_*:=
(\gamma^\Stack_{u',v})^{-1}\odot\St(\beta)_*\odot\gamma^\Stack_{v',u}:
\St(u)_*\circ\St(v')_*\Rightarrow\St(v)_*\St(u')_*
$$
for the external square subdiagram. In light of proposition
\ref{prop_Brexit}, we then see that the resulting cubical
diagram commutes on $2$-cells, in the sense of remark
\ref{rem_transit-base-change}(iii).

\sset\subsubsection{}\label{subsec_return-to-wlinks}
In the situation of \eqref{subsec_rotate-the-cube}, we
associate with $v$ the links $(\tilde v,\sV\tdu\bCat)$,
$\sV\tdu\St(v)$ and $\sV\tdu\Fib(v)$ of the $2$-category
$\sV\tdu\overline{2\tdu\bCat}$, as in remark
\ref{rem_get-link-if-u-strong}, as well as $1$-cells
$$
((-)^a_D,(-)^a_C,\Upsilon(\cD^v)):\sV\tdu\Fib(v)\to\sV\tdu\St(v)
\qquad
(\cFib_\cD,\cFib_\cC,\perp^v):
(\tilde v,\sV\tdu\bCat)\to\sV\tdu\Fib(v)
$$
in the $2$-category $\wLink(\sV\tdu\overline{2\tdu\bCat})$, whose
composition is the $1$-cell :
$$
(\cFib_D^a,\cFib^a_C,\Delta^v):
(\tilde v,\sV\tdu\bCat)\to\sV\tdu\St(v)
$$
(and likewise for $v'$). Additionally, we get two more $1$-cells
of $\wLink(\sV\tdu\overline{2\tdu\bCat})$ :
$$
\begin{aligned}
((\tilde u{}',\bCat)_*,(\tilde u,\bCat)_*,(\tilde\beta,\bCat)_*):&\,
(\tilde v{}',\bCat)\to(\tilde v,\bCat) \\
(\St(u')_*,\St(u)_*,\St(\beta)_*^\gamma):&\,\St(v')\to\St(v).
\end{aligned}
$$

\begin{corollary}\label{cor_conditions-for-trivial-bc}
With the notation of \eqref{subsec_return-to-wlinks}, let
$\cA_\bullet\in\Ob((C',\bCat)^\sim)$. We have :
\begin{enumerate}
\item
If\/ $\cFib_{\cC'}(\cA_\bullet)$ is a stack on $C'$, then
$\Delta^u_{\cA_\bullet}$ is an equivalence of categories.
\item
If\/ $u$ is also cocontinuous, then $\Delta^u$ is a pseudo-natural
equivalence.
\item
$\Upsilon(\cFib^a_D,\cFib^a_C,\Delta^v)$ is a
pseudo-natural equivalence.
\item
In the situation of \eqref{subsec_rotate-the-cube}, suppose
that both $\cFib_{\cC'}(\cA_\bullet)$ and
$\cFib_{\cD'}((\tilde v{}',\bCat)^*\cA_\bullet)$ are stacks.
Then $\cFib^a_D(\Upsilon((\tilde\beta,\bCat)_*)_{\cA_\bullet})$
is an equivalence if and only if the same holds for
$\Upsilon(\St(\beta)^\gamma_*)_{\cFib(\cA_\bullet)}$.
\end{enumerate}
\end{corollary}
\begin{proof}(i): It suffices to show that $\Upsilon(\cD^u)_\cE$
is an equivalence if $\cE$ is a stack on $C'$. Thus, let
$(\eta^C,\eps^C)$ be the unit and counit for the $2$-adjoint
pair $(()^a_C,i_C)$ (notation of
\eqref{subsec_weak-morph-of-sites-and-cocont}) and define
likewise $(\eta^{C'},\eps^{C'})$. Since $i_C$ is fully faithful,
$\eps^C$ is a pseudo-natural equivalence, and likewise for
$\eps^{C'}$. From the triangular identities of theorem
\ref{th_2-adjunction}(i) we deduce that also $\eta^C*i_C$
and $\eta^{C'}*i_{C'}$ are pseudo-natural equivalences; the
assertion follows directly.

(ii): Again, it suffices to check that $\Upsilon(\cD^u)$
is a pseudo-natural equivalence if $u$ is cocontinuous;
this is known by corollary \ref{cor_weak-morph-of-sites-and-cocont}.

\begin{claim}\label{cl-reduce-to-small-lex}
In order to prove (iii), we may assume that $v$ is a morphism
of small lex-sites.
\end{claim}
\begin{pfclaim} Let $T$ (resp. $T'$) be the site whose
underlying category is $C^\sim_\sU$ (resp $D^\sim_\sU$),
with its canonical topology. Pick a universe $\sV'$ with
$\sV\subset\sV'$, such that $T$ and $T'$ are $\sV'$-small.
By direct inspection, we see that $\sV\tdu\Delta^v$ is the
restriction of $\sV'\tdu\Delta^v$; hence if
$\Upsilon(\sV'\tdu\Delta^v)$ is a pseudo-natural equivalence,
the same holds for $\Upsilon(\sV\tdu\Delta^v)$. We may then replace
$\sV$ by $\sV'$, and suppose that $C,D,T,T'$ are $\sV$-small;
notice that $w:=\tilde u{}^*:T'\to T$ is a morphism of lex-sites
(remark \ref{rem_choose-two-univs}(ii)). We consider the induced
essentially commutative diagram of morphisms of sites :
$$
\xymatrix{T' \ar[r]^-w \ar[d]_{h^a_\cD}
\drtwocell\omit{_\ \alpha} &
T \ar[d]^{h^a_\cC} \\
D \ar[r]^-v  & C
}$$
(lemma \ref{lem_cont-funct-site}(ii)) which, according to
\eqref{subsec_rotate-the-cube}, induces a diagram :
$$
\xymatrix@C+20pt{
(D,\sV\tdu\bCat)^\sim \ar[rrr]^-{(\tilde v,\sV\tdu\bCat)_*}
\ar[ddd]_{\cFib^a_D} & &
\dltwocell\omit{_\ \ \ \ \ \ \ \ \ \ \ \ \ \ \ \ \ \ (\tilde\alpha,\bCat)_*}
& (C,\sV\tdu\bCat)^\sim \ar[ddd]^{\cFib^a_C} \\
&  \dltwocell\omit{^\ \Delta^{h^a_\cD}\ \ \ \ \ \ \ \ \ \ \ \ \ \ }
(T',\sV\tdu\bCat)^\sim \ar[r]^-{(\tilde w,\sV\tdu\bCat)_*}
\ar[d]_{\cFib^a_{T'}} \drtwocell\omit{_\ \ \ \Delta^w}
\ar[lu]_-{(\tilde h^a_\cD,\bCat)_*} & (T,\sV\tdu\bCat)^\sim
\ar[d]^{\cFib^a_T} \ar[ru]^-{(\tilde h^a_\cC,\bCat)_*}
\drtwocell\omit{_\ \ \ \ \ \Delta^{h^a_\cC}} \\
& \sV\tdu\Stack(T') \ar[r]_-{\sV\tdu\St(w)_*}
\ar[ld]_-{\sV\tdu\St(h^a_\cD)_*}
\drtwocell\omit{_\ \St(\alpha)^\gamma_*\ \ \ \ \ \ \ \ \ \ \ \ \ \ \ } &
\sV\tdu\Stack(T) \ar[rd]^-{\sV\tdu\St(h^a_\cC)_*} & \\
\sV\tdu\Stack(D) \ar[rrr]_-{\sV\tdu\St(v)_*} & & &
\sV\tdu\Stack(C)
}$$
completed by adding the orientation $\Delta^v$ for the
external square subdiagram. Indeed, this diagram is obtained
by rotating suitably the corresponding diagram of
\eqref{subsec_rotate-the-cube}, and notice that, after
such rotation, the diagram still commutes on $2$-cells,
but now in the sense of \eqref{subsec_transfer-base-change}.
The four diagonal arrows of the diagrams are
$2$-equivalences, by theorem \ref{th_canon-topos}(iii),
proposition \ref{prop_equiv-and-2-equiv}(ii), and remark
\ref{rem_alternative-inverse-image}(ii). Then
$\Upsilon((\tilde\alpha,\bCat)_*)$ and $\Upsilon(\St(\alpha)^\gamma_*)$
are pseudo-natural equivalences, by remark
\ref{rem_when-Upsilon-inverts}. Combining with (ii),
theorem \ref{th_canon-topos}(iii) and remark
\ref{rem_transit-base-change}(i), we conclude that if
$\Upsilon(\Delta^w)$ is a pseudo-natural equivalence,
the same holds for $\Upsilon(\Delta^v)$, whence the claim.
\end{pfclaim}

(iii): By claim \ref{cl-reduce-to-small-lex}, we shall
henceforth suppose that $v$ is a morphism of small lex-sites.
By remark \ref{rem_double-upsilon} we have
$\Upsilon(\Upsilon(\cD^v))=\one_{i_C\circ\sV\tdu\St(v)_*}^\dagger$,
and since
$$
\one_{i_C\circ\sV\tdu\St(v)_*}:
\Fib(v)\circ((-)^a_{C'},i_{C'},\eta^C,\eps^C)\Rightarrow
((-)^a_C,i_C,\eta^C,\eps^C)\circ\St(v)
$$
is an invertible $2$-cell in the category
$\sLink(\sV\tdu\overline{2\tdu\bCat})$, the strict isomorphisms
of proposition \ref{prop_isoms-of-links} show that
$\one_{i_C\circ\sV\tdu\St(v)_*}^\dagger$ is invertible as well.
On the other hand, $\Upsilon(\perp^v)$ is not necessarily a
pseudo-natural equivalence, but in order to conclude the
proof of (iii) it will suffice to show that
$(-)^a_D*\Upsilon(\perp^v)$ is a pseudo-natural equivalence.
To this aim, notice that \eqref{eq_with-perps} (for the morphism
$v$) can be further decomposed as a chain of three
oriented squares :

$$
\xymatrix@C+60pt{
(D,\sV\tdu\bCat)^\sim \ar[r]^-{(\tilde v,\sV\tdu\bCat)_*}
\ar[d]_{j_D} \drtwocell\omit{_\ \ \ \perp^v_1} &
(C,\sV\tdu\bCat)^\sim \ar[d]^{j_C} \\
(\cD,\sV\tdu\bCat)^\wedge \ar[r]^-{(v,\sV\tdu\bCat)^\wedge}
\ar[d]_{j_\cD} \drtwocell\omit{_\ \ \ \perp^v_2} &
(\cC,\sV\tdu\bCat)^\wedge \ar[d]^{j_\cC} \\
\sPsFun(\cD^o,\sV\tdu\bCat) \ar[r]^-{\sPsFun(v^o,\sV\tdu\bCat)}
\ar[d]_{\cFib_\cD} \drtwocell\omit{_\ \ \ \perp^v_3} &
\sPsFun(\cC^o,\sV\tdu\bCat) \ar[d]^{\cFib_\cC} \\
\sV\tdu\Fib(\cD) \ar[r]_-{\sV\tdu\Fib(v)^*} & \sV\tdu\Fib(\cC)
}$$
where $j_C$, $j_\cC$, $j_D$ and $j_\cD$ are the inclusion
pseudo-functors, and where now $\cFib_\cC$ and $\cFib_\cD$
are strong $2$-equivalences (theorem
\ref{th_fundamental-fibrations}).
The orientation $\perp^v_3$ is still a pseudo-natural
equivalence, and both orientations $\perp^v_1$ and
$\perp^v_2$ are identities. After choosing as usual a
unit and counit $(\eta,\eps)$ for the $2$-adjoint pair
$((v,\sV\tdu\bCat)_!,(v,\sV\tdu\bCat)^\wedge)$ and a
unit and counit $(\eta',\eps')$ for the $2$-adjoint pair
$(2\tdu\!\!\int^{v^o},\sPsFun(v^o,\sV\tdu\bCat))$, we may
regard the three squares as oriented squares of weak links.
Then $\Upsilon(\perp^v_3)$ is well defined, and is a
pseudo-natural equivalence (remark
\ref{rem_when-Upsilon-inverts}).  Next, let $\cA_\bullet$
be any presheaf of categories on $\cC$, and set
$\cB_\bullet:=(v,\sV\tdu\bCat)^\wedge\cA$; by definition
we have
$$
\Upsilon(\perp^v_2)_{\cA_\bullet}=\eps'_{\cB_\bullet}\odot
2\tdu\!\!\int^{v^o}\!\!\!\eta_{\cA_\bullet}.
$$
But the discussion of remark \ref{rem_alternative-inverse-image}(v)
shows that
$$
\eps'_{\cB_\bullet}=\eps_{\cB_\bullet}
\qquad\text{and}\qquad
2\tdu\!\!\int^{v^o}\!\!\!\eta_{\cA_\bullet}=
(v,\sV\tdu\bCat)_!(\eta_{\cA_\bullet}).
$$
Hence, the triangular identities of \eqref{subsec_adj-pair}
yield $\Upsilon(\perp^v_2)_{\cA_\bullet}=\one_{\cB_\bullet}$.
Thus, we are further reduced to checking that
$\cFib^a_D*\Upsilon(\perp^v_1)$ is a pseudo-natural equivalence,
or equivalently, that the same holds for
${}^o(\cFib^a_D)*{}^o\Upsilon(\perp^v_1)$. However, let
$[-]^a_C$ (resp. $[-]^a_D$) be the left $2$-adjoint of $j_C$
(resp. of $j_D$), and $(\eta^C,\eps^C)$ (resp. $(\eta^D,\eps^D)$)
the unit and counit of the $2$-adjoint pair $([-]^a_C,j_C)$
(resp. $([-]^a_D,j_D)$); by proposition
\ref{prop_opp-links-and-base-ch} we have
$$
{}^o\Upsilon(\perp^v_1)=(\Upsilon({}^o(\perp^v_1)^\dagger)=
{}^o((j_D\circ(\tilde v,\sV\tdu\bCat)^*)*\eps^C)\odot
(j_D*\perp^{v\dagger}_1*j_C)\odot(\eta^D*((v,\sV\tdu\bCat)_!\circ j_C)))
$$
and notice that $\eps^C$ is a pseudo-natural equivalence, since
$j_C$ is fully faithful (corollary \ref{cor_fully-faith-2-adjoint}).
Also $\perp^{v\dagger}_1$ is a pseudo-natural equivalence, by
virtue of proposition \ref{prop_isoms-of-links}. Therefore,
let  $\cB_\bullet$ be any sheaf of categories on $C$, and set
$\cB'_\bullet:=(v,\sV\tdu\bCat)_!\cB_\bullet$; we are reduced to
checking that the cartesian functor $\cFib^a_D(\eta^C_{\cB'_\bullet})$
is an equivalence of categories. But remark \ref{rem_realign}
easily implies that $\cFib_\cD(\eta^C_{\cB'_\bullet})$ is
$t$-covering for $t=0,1,2$, so the assertion follows from
propositions \ref{prop_unit-sep-cov-faith}(iii) and
\ref{prop_faith-separ-cover}.

(iv): Under the stated assumptions, (i) says that both
$\Delta^u_{\cA_\bullet}$ and $\Delta^{u'}_{(\tilde v{}',\bCat)^*\cA}$
are equivalences, and both $\Upsilon(\Delta^v)$ and
$\Upsilon(\Delta^{v'})$ are pseudo-natural equivalences, by
virtue of (iii). Then the assertion follows by arguing as
in remark \ref{rem_transit-base-change}(i).
\end{proof}

\begin{example}
(i)\ \
Let $C:=(\cC,J)$ be a $\sU$-site, and $A$ any small category.
The {\em constant presheaf of categories on $\cC$ with value
$A$} is the presheaf $A_\cC$ such that $A_\cC(X):=A$ for every
$X\in\Ob(\cC)$, and $A_\cC(f):=\one_A$ for every morphism $f$
of $\cC$. The associated sheaf $A^a_\cC$ is then the {\em
constant sheaf of categories on $C$ with value $A$}. Likewise,
we say that $\cFib(A_\cC)$ is the {\em constant fibration on
$\cC$ with value $A$}, and $\cFib(A_\cC)^a$ is the {\em constant
stack on $C$ with value $A$}.

(ii)\ \
Now, let $C':=(C',J')$ be another $\sU$-site, and $u:C'\to C$
a morphism of sites; we notice that $(\tilde u,\bCat)^*A_\cC$
is isomorphic to the constant sheaf of categories on $C'$
with value $A$. Indeed, for the proof we are easily reduced
to the corresponding assertion for constant sheaves of sets,
which is known by remark \ref{rem_morph-of-sites}(iii).

(iii)\ \
In light of (i) and lemma \ref{lem_stacks-and-sheaves-of-cats}(iii),
we deduce that the stack $\St(u)^*(\cFib(A_\cC)^a)$ is equivalent to
the constant stack on $C'$ with value $A$.
\end{example}

\sset\subsubsection{}\label{subsec_from-beta-dagger-to-tilde}
In section \eqref{sec_C2-for-fibred-sites} we will prove two
{\em base change theorems}, asserting that, for suitable oriented
diagrams of sites as in \eqref{subsec_rotate-the-cube}, the
{\em base change map} $\Upsilon(\St(\beta)^\gamma_*)$ is a
pseudo-natural equivalence. To this aim we shall apply corollary
\ref{cor_conditions-for-trivial-bc}(iv), thereby reducing the
assertion to the corresponding one for
$\Upsilon((\tilde\beta,\bCat)_*)$, {\em i.e.} we shall deduce
a base change theorem for stacks from one for sheaves of
categories. The latter in turn can be reduced to the
corresponding assertion for $\Upsilon(\beta^\sim_*)$, as we explain
hereafter. Indeed, by remark \ref{rem_Ob-and-Morph}(v,vi,vii)
we deduce from \eqref{subsec_rotate-the-cube} a diagram of
$2$-categories whose four diagonal arrows are $2$-equivalences :
$$
\xymatrix{ (D',\bCat)^\sim \ar[rrrrr]^-{(\tilde v',\bCat)_*}
\ar[ddd]_{(\tilde u',\bCat)_*} & & & & &
(C',\bCat)^\sim \ar[ddd]^{(\tilde u,\bCat)_*} \\
& \bCat^*(D'^\sim) \ar[rrr]^-{\bCat^*(\tilde v{}'_*)} \ar[lu]
\ar[d]_{\bCat^*(\tilde u{}'_*)} &
\drtwocell\omit{_\ \ \ \ \ \ \ \ \ \ \bCat^*(\tilde\beta_*)} & &
\bCat^*(C'^\sim) \ar[d]^{\bCat^*(\tilde u_*)} \ar[ru] \\
& \bCat^*(D^\sim) \ar[rrr]_-{\bCat^*(\tilde v_*)} \ar[dl] & & &
\bCat^*(C^\sim) \ar[dr] \\
(D,\bCat)^\sim \ar[rrrrr]^-{(\tilde v,\bCat)_*}
& & & & & (C,\bCat)^\sim
}$$
whose four external trapezoidal subdiagrams are strictly
commutative, and can therefore be oriented by adding in
the respective identity pseudo-natural transformations.
We further complete the diagram by inserting the orientation
$$
(\tilde\beta,\bCat)_*:
(\tilde u,\bCat)_*\circ(\tilde v{}',\bCat)_*
\Rightarrow(\tilde v,\bCat)_*\circ(\tilde u{}',\bCat)_*
$$
for the ``front face'' of the resulting cubical diagram $\cD$.
We regard $\cD$ as a diagram of $1$-cells and $2$-cells in
$\sV\tdu\overline{2\tdu\bCat}$ (for a suitable universe $\sV$);
then it is clear that $\cD$ commutes on $2$-cells in the sense
of \eqref{subsec_transfer-base-change}. Moreover, all the
horizontal arrows in $\cD$ admit left $2$-adjoint pseudo-functors,
which we regard as left adjoint $1$-cells, in the $2$-category
$\sV\tdu\overline{2\tdu\bCat}$; after fixing such a system of
left adjoints, we may regard $\cD$ as a diagram oriented squares
of links in $\sV\tdu\overline{2\tdu\bCat}$, and especially, the
$2$-cells $\Upsilon((\tilde\beta,\bCat)_*)$ and
$\Upsilon(\bCat^*(\tilde\beta_*))$ are then well defined.
Invoking again remarks \ref{rem_when-Upsilon-inverts}), and
\ref{rem_transit-base-change}(i) we deduce that
$\Upsilon((\tilde\beta,\bCat)_*)$ is invertible if and
only if the same holds for $\Upsilon(\bCat^*(\tilde\beta_*))$
(recall that the invertible $2$-cells of
$\sV\tdu\overline{2\tdu\bCat}$ are the equivalence classes of
pseudo-natural equivalences of pseudo-functors). According
to proposition \ref{prop_functoriality-of-Link}, we have 
$\Upsilon(\bCat^*(\tilde\beta_*))=\bCat^*(\Upsilon(\beta_*^\sim))$,
where $\Upsilon(\beta_*^\sim)$ is the base change transformation
associated with the oriented square of links
\eqref{eq_links-from-topoi}; so, if $\Upsilon(\beta^\sim_*)$ is
invertible, the same holds for $\Upsilon((\tilde\beta,\bCat)_*)$.

\sset\subsubsection{}\label{subsec_pashimotta}
Here is a first illustration of the method explained in
\eqref{subsec_from-beta-dagger-to-tilde}. Consider a
morphism of sites $u:C':=(\cC',J')\to C:=(\cC,J)$; we
assume that all finite products are representable in $\cC$
and $\cC'$, and that $u$ commutes with such products. For
every $X,Y\in\Ob(\cC)$ we choose a representative
$X\times Y\in\Ob(\cC)$ for the product of $X$ and $Y$,
and let $X\xleftarrow{q_{X,Y}}X\times Y\xrightarrow{p_{X,Y}}Y$
be the universal projections. For every $X\in\Ob(\cC)$ we get
a morphism of sites $p_X:C/X\to C$ that assigns to every $Y\in\Ob(\cC)$
the projection $p_{X,Y}$ : see remark \ref{rem_continue-local}(iii).
Likewise we define the morphism of sites $p_{uX}:C'/uX\to C'$.
For every $Y\in\Ob(\cC)$, there is by assumption a unique
isomorphism $\beta_Y:u(X\times Y)\isom uX\times uY$ in $\cC'$
such that $p_{uX,uY}\circ\beta_Y=u(p_{X,Y})$, and
$q_{uX,uY}\circ\beta_Y=u(q_{X,Y})$, where
$uX\xleftarrow{q_{uX,uY}}uX\times uY\xrightarrow{p_{uX,uY}}uY$
are likewise the universal projections. Especially,
$\beta_Y/X:u_{|X}(p_{X,Y})\isom p_{uX,uY}$ is an isomorphism
in $\cC'/uX$ (notation of \eqref{eq_restrict-over-X}). Thus,
we get an oriented square diagram of sites :
$$
\xymatrix{ C'/uX \ar[r]^-{p_{uX}} \ar[d]_{u_{|X}}
\drtwocell\omit{_\ \beta} & C' \ar[d]^u \\
C/X \ar[r]^-{p_X} & C.
}$$

\begin{proposition}\label{prop_pashimotta}
In the situation of \eqref{subsec_pashimotta}, the base
change transformation $\Upsilon(\St(\beta)^\gamma_*)$ is
a pseudo-natural equivalence.
\end{proposition}
\begin{proof} To begin with, we remark :

\begin{claim}\label{cl_Roger-Moore-died-today}
Let $\cA_\bullet\in\Ob((C',\bCat)^\sim)$ such that
$\cE:=\cFib_{\cC'}(\cA_\bullet)$ is a stack on $C'$.
Then the fibration
$\cFib_{\cC'/uX}(\tilde p_{uX},\bCat)^*\cA_\bullet$
is a stack on $C'/uX$.
\end{claim}
\begin{pfclaim} The isomorphism of functors
$\tilde p{}^*_{uX}\simeq\tilde\ss_{uX*}$ of remark
\ref{rem_continue-local}(iii) induces an isomorphism of
pseudo-functors $(\tilde p_{uX},\bCat)^*\isom(\tilde\ss_{uX},\bCat)_*$,
hence it suffices to show that
$\cF:=\cFib_{\cC'/uX}(\tilde\ss_{uX},\bCat)_*\cA_\bullet=
\cFib_{\cC'/uX}(\cA_\bullet\circ\ss^o_{uX})\simeq\cC'/uX\times_{\cC'}\cE$
is a stack. Thus, let $(Z\xrightarrow{h}uX)\in\Ob(\cC'/uX)$ and
$\cS'\subset(\cC'/uX)/h$ a sieve covering $h$ for the site $C'/uX$;
under the natural identification of categories
$$
(\cC'/uX)/h\isom\cC'/Z
$$
the sieve $\cS'$ corresponds to a sieve $\cS\subset\cC'/Z$
covering $Z$ for the site $C$ (see \eqref{sec_Localization-topoi}).
Then the restriction functor
$$
\cF(h):=\sCart_{\cC'/uX}((\cC'/uX)/h,\cF)\to
\sCart_{\cC'/uX}(\cS',\cF)
$$
is naturally identified with the restriction functor
$$
\cE(Z):=\sCart_{\cC'}(\cC'/Z,\cE)\to\sCart_{\cC'}(\cS,\cE)
$$
whence the contention.
\end{pfclaim}

By claim \ref{cl_Roger-Moore-died-today}, lemma
\ref{lem_stacks-and-sheaves-of-cats}(ii), corollary
\ref{cor_conditions-for-trivial-bc}(iv) and the discussion of
\eqref{subsec_from-beta-dagger-to-tilde}, we are reduced to
showing that the base change transformation
$\Upsilon(\beta^\sim_*):\tilde p{}^*_X\circ\tilde u_*\to
\tilde u_{|X*}\circ\tilde p{}^*_{uX}$ is an isomorphism of functors.
Now, recall that the source functor $\ss_X:\cC/X\to\cC$ is
continuous for the sites $C$ and $C/X$, and is left adjoint
to $p_X$; more precisely, we have an explicit adjunction :
$$
\theta_{h,Z}:\Hom_\cC(Y,Z)\isom\Hom_{\cC/X}(h,p_{X,Z})
\quad
\text{for every $(Y\xrightarrow{h}X)\in\Ob(\cC/X)$ and $Z\in\Ob(\cC)$}
$$
that assigns to every morphism $f:Y\to Z$ of $\cC$ the unique
morphism $f^*/X:h\to p_{X,Z}$ of $\cC/X$ such that $q_{X,Z}\circ f^*=f$.
The adjunction $\theta_{h,Z}$ induces a corresponding adjunction
for the pair $(\tilde\ss_{X*},\tilde p_{X*})$, as described in
remark \ref{rem_opposite-Fun}(iii,iv). Explicitly, the unit
of this induced adjunction assigns to every sheaf $\cF$ on
$C$ the morphism of sheaves :
$$
\eta^X_\cF:\cF\to\cF\circ\ss^o_X\circ p^o_X
\qquad
Y\mapsto(\cF(q_{X,Y}):\cF(Y)\to\cF(X\times Y))
$$
and the counit assigns to every sheaf $\cG$ on $C/X$ the
morphism of sheaves :
$$
\eps^X_\cG:\cG\circ p^o_X\circ\ss^o_X\to\cG
\qquad
(Y\xrightarrow{h}X)\mapsto(\cG(h^*/X):\cG(p_{X,Y})\to\cG(h)).
$$
The same description applies to the adjunction for the pair
$(\tilde\ss_{uX*},\tilde p_{uX*})$, and we get therefore
isomorphisms of links
$$
\begin{aligned}
\cL_X:=&(\tilde\ss_{X*},\tilde p_{X*},\eta^X,\eps^X)\isom
\tilde p_X:=(\tilde p{}^*_X,\tilde p_{X*},\eta^{p_X},\eps^{p_X}) \\
\cL_{uX}:=&(\tilde\ss_{uX*},\tilde p_{uX*},\eta^{uX},\eps^{uX})\isom
\tilde p_{uX}:=(\tilde p{}^*_{uX},\tilde p_{uX*},\eta^{p_{uX}},\eps^{p_{uX}}).
\end{aligned}
$$
Hence, it suffices to check that the base change transformation
for the oriented diagram of links:
$$
\cD\qquad :\qquad
{\spreaddiagramcolumns{+20pt}\diagram
(C'/uX)^\sim \ar[r]^-{\cL_{uX}} \ar[d]_{\tilde u_{|X}}
\drtwocell\omit{_\ \ \ \tilde\beta_*} & C'^\sim \ar[d]^{\tilde u} \\
(C/X)^\sim \ar[r]_-{\cL_X} & C^\sim
\enddiagram}\qquad\qquad
$$
is an isomorphism of functors $\Upsilon(\cD):
\tilde\ss_{X*}\circ\tilde u_*\isom\tilde u_{|X*}\circ\tilde\ss_{uX*}$.
However, a direct computation that we leave to the reader easily
shows that $\Upsilon(\cD)$ is the identity automorphism of
$\tilde\ss_{X*}\circ\tilde u_*=\tilde u_{|X*}\circ\tilde\ss_{uX*}$
\end{proof}

\subsection{Stacks in groupoids and ind-finite stacks}
Let $C:=(\cC,J)$ be any site; for every universe $\sV$ denote by
$$
\sV\tdu\Stack^\times(C)
\qquad\text{and}\qquad
\sV\tdu\StGpd(C)
$$
respectively the strong $2$-subcategory of
$\sV\tdu\Fib^\times(\cC)$ whose objects are the $\sV$-stacks
on $C$ and the strong $2$-subcategory of $\sV\tdu\Stack(C)$
whose objects are the {\em $\sV$-stacks in groupoids} on $C$,
{\em i.e.} the $\sV$-stacks on $C$ that are fibrations in
groupoids. For every $\cC$-fibration $F:\cE\to\cC$, every
$X\in\Ob(\cC)$, and every sieve $\cS\subset\cC/X$, the
inclusion $\cE^\times\to\cE$ yields a commutative diagram
$$
\xymatrix{
\cE(X)^\times \ar[r] \ar[d] & \cE^\times(X) \ar[d] \\
\sCart_\cC(\cS,\cE)^\times \ar[r] &
\sCart_\cC(\cS,\cE^\times)
}$$
whose horizontal arrows are isomorphisms of categories,
due to remark \ref{rem_groupoids}(i); notice that the
source functor $\ss:\cS\to\cC$ is a fibration in groupoids,
by remark \ref{rem_sieves-and-sub}(ii). Suppose now that
$\cE$ is $i$-separated on $C$ for some $i\in\{0,1,2\}$;
from remark \ref{rem_cat-groupoid}(ii) we deduce that the
left vertical arrow is an $i$-faithful functor whenever
$\cS\in J(X)$, and then the same holds for the right
vertical arrow, so finally $\cE^\times$ is also $i$-separated
on $C$. Especially, the pseudo-functor $(-)^\times_\cC$
restricts to a pseudo-functor
$$
(-)^\times_C:\sV\tdu\Stack^\times(C)\to\sV\tdu\StGpd(C).
$$
Next, if $C$ is a small site, taking into account remark
\ref{rem_cat-groupoid}(iii) we get a natural identification:
$$
(\cE^+_X)^\times\isom\colim_{\cS\in J(X)}\sCart_\cC(\cS,\cE)^\times
\isom\colim_{\cS\in J(X)}\sCart_\cC(\cS,\cE^\times)=((\cE^\times)^+)_X
$$
whence natural isomorphisms of $\cC$-fibrations :
\set\begin{equation}\label{eq_gazouille}
(\cE^+)^\times\isom(\cE^\times)^+
\qquad
(\cE^a)^\times\isom(\cE^\times)^a.
\end{equation}
If $C':=(\cC',J')$ is another site, and $u:C\to C'$
a weak morphism of sites, it is also clear that
$\sV\tdu\St(u)_*$ restricts to pseudo-functors
$$
\sV\tdu\St^\times(u)_*:\sV\tdu\Stack^\times(C)\to\sV\tdu\Stack^\times(C')
\qquad
\sV\tdu\StGpd(u)_*:\sV\tdu\StGpd(C)\to\sV\tdu\StGpd(C')
$$
that make commute the diagram of $2$-categories :
\set\begin{equation}\label{eq_burkini}
{\spreaddiagramcolumns{+40pt}\diagram
\sV\tdu\Stack^\times(C) \ar[d]_{(-)^\times_C}
\ar[r]^-{\sV\tdu\St^\times(u)_*} &
\sV\tdu\Stack^\times(C') \ar[d]^{(-)^\times_{C'}} \\
\sV\tdu\StGpd(C) \ar[r]^-{\sV\tdu\StGpd(u)_*} &
\sV\tdu\StGpd(C').
\enddiagram}
\end{equation}

\sset\subsubsection{}\label{subsec_StGpd}
Consider any $\sU$-site $C:=(\cC,J)$, and let $j:\cT\to\cC$
be the inclusion functor of a small subcategory whose set of
objects is a small topologically generating family for $C$;
we endow as usual $\cT$ with the topology $J_\cT$ induced by
$J$ via $j$, so that $j$ is a continuous functor for $J$ and
$J_\cT$. Say that $\cC$ is $\sU'$-small for some universe
$\sU'$, and recall that the rule :
$\cA\mapsto\cA^a:=\sU'\tdu\Fib(j)_*(\sU'\tdu\Fib(j)^*(\cA)^a)$
for every fibration $\cA$ with small fibres on $\cC$ yields
a left $2$-adjoint for the inclusion pseudo--functor
$\Stack(C)\to\Fib(\cC)$ (see the proof of corollary
\ref{cor_stackificate-U-sites}). Combining with
\eqref{subsec_times-and-upper-*} and \eqref{subsec_times-and-lower-*}
we deduce a pseudo-commutative diagram :
$$
\xymatrix{ \sV\tdu\Fib^\times(\cC) \ar[r]^-{(-)^a}
\ar[d]_{(-)^\times_\cC} & \sV\tdu\Stack^\times(C) \ar[d]^{(-)^\times_C} \\
\sV\tdu\Gpd(\cC) \ar[r]^-{(-)^a} & \sV\tdu\StGpd(C)
}$$
for every universe $\sV$ containing $\sU$.

\begin{proposition}\label{prop_Van-Gogh}
For every morphism of $\sU$-sites $u:C':=(\cC',J')\to C:=(\cC,J)$,
and every universe $\sV$ containing $\sU$, the pseudo-functor
$\sV\tdu\St(u)^*$ restricts to pseudo-functors
$$
\sV\tdu\St^\times(u)^*:\sV\tdu\Stack^\times(C)\to\sV\tdu\Stack^\times(C')
\qquad
\sV\tdu\StGpd(u)^*:\sV\tdu\StGpd(C)\to\sV\tdu\StGpd(C')
$$
and we have a pseudo-commutative diagram :
$$
\xymatrix@C+40pt{
\sV\tdu\Stack^\times(C) \ar[r]^-{\sV\tdu\St^\times(u)^*}
\ar[d]_{(-)^\times_C} &
\sV\tdu\Stack^\times(C') \ar[d]^{(-)^\times_{C'}} \\
\sV\tdu\StGpd(C) \ar[r]^-{\sV\tdu\StGpd(u)^*} & \sV\tdu\StGpd(C').
}$$
\end{proposition}
\begin{proof} Endow the topoi $T:=C^\sim_\sU$ and $T':=C'^\sim_\sU$
with their canonical topologies; we get an essentially commutative
diagram of morphisms of sites :
$$
{\diagram
T' \ar[r]^-{u'} \ar[d]_{h^a_{C'}} &
T \ar[d]^{h^a_C} \\
C' \ar[r]^-u & C
\enddiagram}\qquad
\text{with $u':=\tilde u{}^*$}
$$
(lemma \ref{lem_cont-funct-site}(ii)). After replacing
$\sU$ by a larger universe, we may assume that $T$ and
$T'$ are small; notice also that $u'$ is a morphism
of lex-sites (remark \ref{rem_choose-two-univs}(ii)).
Especially, the proposition is already known for $u'$,
by virtue of \eqref{eq_awdnews} and \eqref{eq_gazouille}.

\begin{claim}\label{cl_pocket}
The proposition holds for $h^a_C$ and $h^a_{C'}$.
\end{claim}
\begin{pfclaim} It suffices to check the assertion for
$h^a_C$, since the same argument will work for $h^a_{C'}$.
Now, let $\cA$ be a stack in groupoids on $C$; we need
to check that $\St(h^a_C)^*(\cA)$ is a stack of groupoids
on $T$. However, since $\St(h^a_C)_*$ is a $2$-equivalence
(proposition \ref{prop_equiv-and-2-equiv}(ii)), we may
assume that $\cA=\St(h^a_C)_*\cB$ for a stack $\cB$ on
$T$, and we may moreover assume that $\cB$ is a stack
in groupoids, by \eqref{eq_gazouille}. But since
$\St(h^a_C)^*(\cA)$ is equivalent to $\cB$, the assertion
follows. This shows that the pseudo-functor $\StGpd(h^a_C)^*$
is well defined, and then it is clearly a pseudo-inverse
for $\StGpd(h^a_C)_*$. But then the required pseudo-commutativity
of the resulting diagram follows easily from the commutativity
of \eqref{eq_burkini} : details left to the reader.
\end{pfclaim}

Now, let $\cA$ be a stack in groupoids on $C$; we need to
check that $\cB:=\St(u)^*(\cA)$ is a stack in groupoids on
$C'$, and by claim \ref{cl_pocket} it suffices to show that
$\St(h^a_{C'})^*(\cB)$ is a stack in groupoids on $T'$; but
the latter is equivalent to $\St(u')^*\circ\St(h^a_C)^*(\cA)$,
and the proposition is already known for both $u'$ and $h^a_C$.
This shows that $\StGpd(u)^*$ is well defined.

Lastly, we consider the diagram of $2$-categories :
$$
\xymatrix@C+40pt{\StGpd(T) \ar[rrr]^-{\StGpd(u')^*}
\ar[ddd]|{\StGpd(h^a_C)^*} & & &
\StGpd(T') \ar[ddd]|{\StGpd(h^a_{C'})^*} \\
& \Stack^\times(T) \ar[r]^-{\St^\times(u')^*}
\ar[lu]^-{(-)^\times_T} \ar[d]_{\St^\times(h^a_C)^*} &
\Stack^\times(T') \ar[d]^{\St^\times(h^a_{C'})^*} \ar[ru]_-{(-)^\times_{T'}} \\
& \Stack^\times(C) \ar[r]^-{\St^\times(u)^*}
\ar[ld]_-{(-)^\times_C} & \Stack^\times(C') \ar[rd]^-{(-)^\times_{C'}} \\
\StGpd(C) \ar[rrr]^-{\StGpd(u)^*} & & & \StGpd(C')
}$$
whose inner and outer square subdiagrams are pseudo-commutative.
Moreover, by the foregoing, we know that also the left, the
right, and the top trapezoidal subdiagrams pseudo-commute.
Since $\St^\times(h^a_C)$ is a $2$-equivalence, a little diagram
chase shows easily that the bottom trapezoidal subdiagram
pseudo-commutes as well, whence the proposition.
\end{proof}

\sset\subsubsection{}\label{subsec_parthenogenesis}
Let $C:=(\cC,J)$ be any $\sU$-site; according to example
\ref{ex_sheaves-as-stacks}, for every universe $\sV$ containing
$\sU$, the pseudo-functor $\cFib_\cC$ of
\eqref{subsec_pi_0-on-groupoids} restricts to a (strict)
pseudo-functor :
$$
\cFib_C:C^\sim_\sV\to\sV\tdu\StGpd(C)
$$
and then it is clear that $\pi_0^\cC$ induces a left $2$-adjoint
pseudo-functor
$$
\pi_0^C:\sV\tdu\StGpd(C)\to C^\sim_\sV
\qquad
\cE\mapsto(\pi_0^\cC(\cE))^a.
$$

\begin{lemma}\label{lem_parthenon}
In the situation of \eqref{subsec_parthenogenesis}, let
$\phi:\cE\to\cF$ be a cartesian functor of fibrations in
groupoids over $\cC$. The following holds :
\begin{enumerate}
\item
$\phi$ is $0$-covering (for the topology of $C$) if and only
if $\pi^\cC_0(\phi)^a:\pi_0^\cC(\cE)^a\to\pi_0^\cC(\cF)^a$ is
an epimorphism of sheaves on $C$.
\item
If $\phi$ is $i$-covering for $i=0,1$, then $\pi^\cC_0(\phi)^a$
is an isomorphism.
\end{enumerate}
\end{lemma}
\begin{proof} Both assertions follow by direct inspection of
the definitions, taking into account corollary \ref{cor_bicover}
and remark \ref{rem_iprippi}(ii,iii).
\end{proof}

\sset\subsubsection{}\label{subsec_pull-back-pi_0^C}
Let $u:C':=(\cC',J')\to C:=(\cC,J)$ be a morphism of
$\sU$-sites; then for every universe $\sV$ containing $\sU$
the induced diagram of $2$-categories
$$
\xymatrix@C+40pt{
\sV\tdu\StGpd(C) \ar[r]^-{\sV\tdu\StGpd(u)^*} \ar[d]_{\pi_0^C} &
\sV\tdu\StGpd(C') \ar[d]^{\pi_0^{C'}} \\
C^\sim \ar[r]^-{\tilde u{}^*} & C'^\sim
}$$
is essentially commutative. Indeed, lemma \ref{lem_parthenon}(ii)
and proposition \ref{prop_unit-sep-cov-faith}(i) imply that for
every fibration $\cE$ over $\cC'$, the unit of adjunction
$\cE\to\cE^a$ induces an isomorphism
$\pi_0^{\cC'}(\cE)^a\isom\pi_0^{C'}(\cE^a)$ of sheaves on $C'$.
The assertion follows easily, taking into account the essential
commutativity of \eqref{eq_lost-at-sea} : the details are left
to the reader.

\sset\subsubsection{Ind-finite stacks}
It is easy to say when a sheaf of groups $G$ on a topological space
$\cT$ is ind-finite, when $\cT$ admits a basis of quasi-compact
open subsets : in which case, one simply asks that the group of
$U$-sections $G(U)$ is a filtered union of finite groups, for every
quasi-compact open subset $U$ of $\cT$. This definition is however
not suitable for more general topological spaces, nor of course
for arbitrary sites. Before we explain a definition that is
appropriate for the general case, let us introduce a notion
of quasi-compactness for such context:

\begin{definition}\label{def_qcoh-obj-site}
Let $(\cC,J)$ be any site, and $X\in\Ob(\cC)$. We say that $X$
is {\em quasi-compact for the topology $J$}, if for every covering
sieve $\cS\in J(X)$ there exists a finite subset $S\subset\cS$
that generates a sieve covering $X$.
\end{definition}

Obviously this definition recovers the standard one, in
the case of the site of open subsets of any topological space.

\begin{definition}\label{def_ind-finite-sheaf}
Let $C:=(\cC,J)$ be a site, and $G$ a sheaf of groups on $C$.
We say that $G$ is {\em ind-finite} if the following holds.
For every $X\in\Ob(\cC)$ and every finite subset
$\Sigma\subset GX$, there exists a covering family
$(f_i:X_i\to X~|~i\in I)$ for the topology $J$ such that
the set $Gf_i(\Sigma):=\{Gf_i(\sigma)~|~\sigma\in\Sigma\}$
generates a finite subgroup of $GX_i$, for every $i\in I$.
\end{definition}

\begin{remark}\label{rem_ind-finite-and-qc}
The first observation is that if $G$ is ind-finite on the
site $(\cC,J)$, then for every quasi-compact $X\in\Ob(\cC)$
the group $GX$ is {\em ind-finite}, {\em i.e.} it is a
filtered union of finite groups. Indeed, let
$\Sigma\subset GX$ be any finite subset; by assumption
we have a covering family $(X_i\to X~|~i\in I)$ such that
$Gf_i(\Sigma)$ generates a finite group $H_i$ for every
$i\in I$, and since $X$ is quasi-compact, we may assume
that $I$ is a finite set. But the natural map
$GX\to\prod_{i\in I}GX_i$ is injective, and maps the subgroup
$H$ generated by $\Sigma$ into the finite group $\prod_{i\in I}H_i$,
so $H$ is finite, whence the claim.
\end{remark}

\sset\subsubsection{}\label{subsec_scampata-bella}
For every group $G$ and every subset $\Sigma\subset G$, let
us write $\La\Sigma\Ra\subset G$ for the subgroup generated
by $\Sigma$. Definition \ref{def_ind-finite-sheaf} prompts
the following construction. Let $C:=(\cC,J)$ be any site, and
$H$ any presheaf of groups on $\cC$. For every $n\in\N$ and
every $X\in\Ob(\cC)$ we set :
$$
(H^n)_\sff(X):=\{(\sigma_1,\dots,\sigma_n)\in H^n(X)~|~
\La\sigma_1,\dots,\sigma_n\Ra\ \text{is a finite group}\}.
$$
Clearly the rule : $X\mapsto(H^n)_\sff(X)$ yields a
sub-presheaf $H^n_\sff$ of the {\em presheaf of sets} $H^n$,
for every $n\in\N$. If $H$ is a sheaf on $C$, $(H^n)_\sff$
is not necessarily a sheaf, hence we consider the sheaf
$$
(H^n)_\slf
$$
defined as the smallest subsheaf of sets of $H^n$ containing
$(H^n)_\sff$. The latter is also (naturally isomorphic to) the
sheaf of sets $(H^n)_\sff^a$ on $C$ associated with the presheaf
$(H^n)_\sff$. Clearly every morphism $\phi:H\to K$ of presheaves
of groups on $\cC$ induces a morphism of presheaves of sets
$$
(\phi^n)_\sff:(H^n)_\sff\to(K^n)_\sff
\qquad
\text{for every $n\in\N$}.
$$
If $H$ and $K$ are sheaves of groups on $C$, we get also
an induced morphism of sheaves of sets
$$
(\phi^n)_\slf:(H^n)_\slf\to(K^n)_\slf
\qquad
\text{for every $n\in\N$}.
$$
Thus, in the notation of definition
\ref{def_sheaves-with-other-values}, we get by these rules
two well defined functors
$$
(-)^n_\sff:(\cC,\sV\tdu\Grp)^\wedge\to\cC_\sV^\wedge
\qquad
(-)^n_\slf:(C,\sV\tdu\Grp)^\sim\to C^\sim_\sV
\qquad
\text{for every $n\in\N$}
$$
for every universe $\sV$.

\begin{remark}\label{rem_CCCP}
Let $C:=(\cC,J)$ be a site, and $H$ a sheaf of groups on $C$.

(i)\ \
Notice that $(H^n)_\sff$ is a separated presheaf for every
$n\in\N$, since it is a sub-presheaf of the sheaf $H^n$; hence
$(H^n)_\slf=(H^n)^+_\sff$ (claim \ref{cl_plus-construction}(ii)).
This means that for every $X\in\Ob(\cC)$, the set
$(H^n)_\slf(X)$ consists of all
$(\sigma_1,\dots,\sigma_n)\in H^n(X)$ for which there exists
a covering family $(f_i:X_i\to X~|~i\in I)$ for the topology
$J$ such that $\La Gf_i(\sigma_1),\dots,Gf_i(\sigma_n)\Ra$ is a
finite subgroup of $GX_i$, for every $i\in I$.

(ii)\ \
We deduce easily from (i) that $H$ is ind-finite if
and only if $H^n=(H^n)_\slf$ for every $n\in\N$.

(iii)\ \
It also follows from (i) that $(H^n)_\sff(X)=(H^n)_\slf(X)$ for
every quasi-compact object $X$ of $\cC$. Indeed, let
$\underline\sigma:=(\sigma_1,\dots,\sigma_n)\in(H^n)_\slf(X)$,
and pick a covering family $(f_i:X_i\to X~|~i\in I)$ verifying
the condition of (i) relative to $\underline\sigma$; since
$X$ is quasi-compact, we may assume that $I$ is a finite set,
and then arguing as in remark \ref{rem_ind-finite-and-qc} we
deduce that $\underline\sigma\in(H^n)_\sff(X)$.
\end{remark}

\begin{lemma}\label{lem_H_f-respects-bicovers}
Let $(\cC,J)$ be a site, and $\phi:H\to K$ a bicovering
morphism of presheaves of groups on $\cC$. Then
$(\phi^n)_\sff:(H^n)_\sff\to(K^n)_\sff$ is bicovering for every
$n\in\N$.
\end{lemma}
\begin{proof} Let $X\in\Ob(\cC)$ and
$\underline\sigma:=(\sigma_1,\dots,\sigma_n)\in(K^n)_\sff(X)$;
according to remark \ref{rem_iprippi}(ii) we have to exhibit
a covering family $(f_i:X_i\to X~|~i\in I)$ such that
$(Kf_i(\sigma_1),\dots,Kf_i(\sigma_n))$ lies in the image of
the map $(\phi^n)_{\sff,X_i}:(H^n)_\sff(X_i)\to(H^n)_\sff(X_i)$
for every $i\in I$. But by assumption, for every $j=1,\dots,n$
there exists a covering sieve $\cS_j\subset\cC/X$ such that
$Kf(\sigma_j)$ lies in the image of the map
$\phi_Y:HY\to KY$ for every $f:Y\to X$ in $\cS_j$. Then the
sieve $\cS:=\cS_1\cap\cdots\cap\cS_n$ still covers $X$ (remark
\ref{rem_topology}(i)), and for every $g:Y\to X$ in $\cS$ we
see that $K^ng(\underline\sigma)=\phi^n_Y(\underline\tau)$
for some $\underline\tau:=(\tau_1,\dots,\tau_n)\in H^n(Y)$.
To conclude, it then suffices to exhibit for every such $g$
and $\tau$ a covering sieve $\cS'\subset\cS/Y$ such that
$H^ng'(\underline\tau)\in(H^n)_\sff(Y')$ for every $g':Y'\to Y$
in $\cS'$. However, by construction we have
$\phi^n_Y(\underline\tau)\in(K^n)_\sff(Y)$; the latter means
that the subgroup
$\La\phi_Y(\tau_1),\dots,\phi_Y(\tau_n)\Ra\subset KY$ is
finite. This in turn means that there exists an integer
$N\geq 2$, such that for every $\underline j\in\{1,\dots,n\}^N$
we may find an integer $M(\underline j)<N$ and
$\underline j'\in\{1,\dots,n\}^{M(\underline j)}$ with
$\prod_{k=1}^N\phi_Y(\tau_{j_k})=
\prod_{k=1}^{M(\underline j)}\phi_Y(\tau_{j'_k})$. By remark
\ref{rem_iprippi}(iii), this implies that there exists
a covering sieve $\cS'_{\underline j}\subset\cC/Y$ such that
$\prod_{k=1}^NHf(\tau_{j_k})=\prod_{k=1}^{M(\underline j)}Hf(\tau_{j'_k})$
for every $f:Y'\to Y$ in $\cS'_{\underline j}$. Then the sieve
$\cS':=\bigcap_{\underline j\in\{1,\dots,n\}^N}\cS'_{\underline j}$
covers $Y$, and we deduce that the sequence
$H^nf(\underline\tau)$ generates a finite subgroup of $H^nY'$,
for every $f:Y'\to Y$ in $\cS'$; {\em i.e.}
$H^nf(\underline\tau)\in(H^n)_\sff(Y')$, as required.

Lastly, let $\underline\sigma,\underline\sigma'\in(H^n)_\sff(X)$
be two sections such that
$\phi^n_X(\underline\sigma)=\phi^n_X(\underline\sigma')$; according
to remark \ref{rem_iprippi}(iii) we need to exhibit a covering
sieve $\cS\subset\cC/X$ such that
$H^nf(\underline\sigma)=H^nf(\underline\sigma')$ for every
$f:Y\to X$ in $\cS$. But by assumption, for $i=1,\dots,n$
there exists a sieve $\cS_i$ such that $Hf(\sigma_i)=Hf(\sigma'_i)$
for every $f$ in $\cS_i$; clearly $\cS:=\cS_1\cap\cdots\cap\cS_n$
will do.
\end{proof}

\sset\subsubsection{}\label{subsec_choo-choo}
Let now $C:=(\cC,J)$ and $C':=(\cC',J')$ be two sites,
$u:\cC'\to\cC$ a functor, and $H$ a presheaf of groups on $\cC$.
Clearly $u^\wedge H$ is a presheaf of groups on $\cC'$, and we have :
$$
u^\wedge((H^n)_\sff)=((u^\wedge H)^n)_\sff
\qquad
\text{for every $n\in\N$}.
$$

\begin{proposition}\label{prop_gattiglio}
In the situation of \eqref{subsec_choo-choo}, suppose that $u$
is cocontinuous for the topologies $J$ and $J'$, and $H$ is a
sheaf of groups on $C$. Then the following holds :
\begin{enumerate}
\item
$\breve u{}^*H$ is a sheaf of groups on $C'$, and
$\breve u{}^*((H^n)_\slf)=((\breve u{}^*H)^n)_\slf$
for every $n\in\N$.
\item
If $H$ is ind-finite, the same holds for $\breve u{}^*H$.
\end{enumerate}
\end{proposition}
\begin{proof}(i): Since the functor $(-)^a$ is exact, it is clear
that the group law $K\times K\to K$ of every presheaf of groups
$K$ on $\cC'$ yields a group law $K^a\times K^a\to K^a$ for the
sheaf $K^a$ : details left to the reader; the first assertion
of the proposition is an immediate consequence.

To check the sought identity, recall that by definition
$\breve u{}^*((H^n)_\slf)=(u^\wedge((H^n)^a_\sff))^a$, and
$$
(u^\wedge((H^n)^a_\sff))^a=(u^\wedge((H^n)_\sff))^a
$$
by lemma \ref{lem_breve} and corollary \ref{cor_bicover}(iii).
Next, $(u^\wedge((H^n)_\sff))^a=(((u^\wedge H)^n)_\sff)^a$,
by \eqref{subsec_choo-choo}. Lastly,
$(((u^\wedge H)^n)_\sff)^a=((((u^\wedge H)^a)^n)_\sff)^a$,
by lemma \ref{lem_H_f-respects-bicovers}, whence the
contention.

(ii) follows immediately from (i) and remark \ref{rem_CCCP}(ii).
\end{proof}

\begin{lemma}\label{lem_slf-and-u_!}
In the situation of \eqref{subsec_choo-choo},
suppose that the category $X/u\cC'$ is cofiltered for every
$X\in\Ob(\cC)$, and let $K$ be a presheaf of groups on $\cC'$.
Then $u_!K$ is a presheaf of groups on $\cC$, and we have
a natural isomorphism of presheaves :
$$
u_!((K^n)_\sff)\isom((u_!K)^n)_\sff
\qquad
\text{for every $n\in\N$}.
$$
\end{lemma}
\begin{proof} Under the stated assumption, the functor $u_!$
is exact ({\em cp}. the proof of corollary \ref{cor_lable}(i)),
so we easily deduce that $u_!K$ is a presheaf of groups, as in
the proof of proposition \ref{prop_gattiglio}(i).

Next, we get a natural morphism of presheaves
$\omega:u_!((K^n)_\sff)\to((u_!K)^n)_\sff$ as follows. First, the
unit of adjunction is a morphism $\eta_K:K\to u^\wedge u_!K$,
which induces a morphism
$(\eta_K^n)_\sff:(K^n)_\sff\to((u^\wedge u_!K)^n)_\sff$. Then, by
\eqref{subsec_choo-choo} we have
$((u^\wedge u_!K)^n)_\sff=u^\wedge(((u_!K)^n)_\sff)$, and the
resulting morphism $(K^n)_\sff\to u^\wedge(((u_!K)^n)_\sff)$
yields by adjunction the sought morphism.

Explicitly, for every $X\in\Ob(\cC)$, every element of
$u_!((K^n)_\sff)(X)$ is the class $[\underline\sigma]$ of
some $\underline\sigma:=(\sigma_1,\dots,\sigma_n)\in K^n(Y)$,
for an object $f:X\to uY$ of $X/u\cC'$, and
$\omega_X([\underline\sigma])$ is the section
$([\sigma_1],\dots,[\sigma_n])\in((u_!K)^n)_\sff(X)$,
where $[\sigma_i]\in u_!K(X)$ is the class of $\sigma_i$,
for $i=1,\dots,n$.

To show the injectivity of $\omega_X$, let
$[\underline\sigma],[\underline\tau]\in u_!((K^n)_\sff)(X)$
such that
$\omega_X([\underline\sigma])=\omega_X([\underline\tau])$;
we may assume that $\underline\sigma,\underline\tau\in K^nY$
for some morphism $f:X\to uY$ of $\cC'$, and the assumption
means that for $i=1,\dots,n$ there exists a morphism
$X/h_i:(g_i:X\to uY_i)\to f$ in $X/u\cC'$ such that
$Kh_i(\sigma_i)=Kh_i(\tau_i)$ in $KY_i$. Since $X/u\cC'$ is
cofiltered, we may then find an object $g':X\to uY'$ of
$X/u\cC'$ and morphisms $X/h'_i:g'\to g_i$ of $X/u\cC'$ with
$(X/h_i)\circ(X/h'_i)=(X/h_j)\circ(X/h'_j)$ for every
$i,j=1,\dots,n$. Recall that $h_i:Y_i\to Y$ and $h'_i:Y'\to Y_i$
are morphisms in $\cC'$, and the foregoing identity means
that $h:=h_i\circ h'_i=h_j\circ h'_j$ for every such $i,j$.
We conclude that $Kh(\sigma_i)=Kh(\tau_i)$ for $i=1,\dots,n$,
whence $[\underline\sigma]=[\underline\tau]$, as required.

To check the surjectivity of $\omega_X$, let
$([\sigma_1],\dots,[\sigma_n])\in((u_!K)^n)_\sff(X)$; hence,
for $i=1,\dots,n$ there exists an object $f_i:X\to uY_i$ of
$X/u\cC'$ such that $\sigma_i\in KY_i$. Arguing as in the
foregoing, we easily reduce to the case where
$Y:=Y_1=\cdots=Y_n$ and $f:=f_1=\cdots=f_n$. For every
morphism $X/h:(g:X\to uY')\to f$ in $X/u\cC'$, let
$G_h:=\La Kh(\sigma_1),\dots,Kh(\sigma_n)\Ra\subset KY'$, and
set $G:=\La[\sigma_1],\dots,[\sigma_n]\Ra\subset u_!K(X)$.
The source morphism $(X/u\cC')/f\to X/u\cC'$ is final
(example \ref{ex_filtered-final}(i)), and a direct inspection
shows that the group $G$ is the colimit of the functor
$$
G_\bullet:((X/u\cC')/f)^o\to\Grp
\qquad
h\mapsto G_h.
$$
To every pair of objects $X/h:(g:X\to uY')\to f$ and
$X/h:(g:X\to uY'')\to f$ of $(X/u\cC')/f$, and every morphism
$(X/l)/f:(X/h')\to(X/h)$ in $(X/u\cC')/f$, the functor
$G_\bullet$ assigns the group homomorphism $G_h\to G_{h'}$
given by the restriction of $Kh:KY'\to KY''$. By assumption,
$G$ is a finite group; let us then remark more generally :

\begin{claim}\label{cl_fg-groups-are-fp}
Let $\Gamma_\bullet:I\to\Grp$ be a functor from a filtered
category $I$, such that :
\begin{enumerate}
\alphaenu
\item
$\Gamma_i$ is a finitely generated group for every $i\in\Ob(I)$
\item
$\Gamma_\phi$ is a surjective group homomorphism for every
morphism $\phi$ of $I$
\item
The colimit $G$ of $\Gamma_\bullet$ is a finite group.
\end{enumerate}
Let also $\tau_\bullet:\Gamma_\bullet\Rightarrow c_G$ be a universal
cocone. Then there exists a cofinal functor $\psi:J\to I$ such
that $\tau*\psi$ is an isomorphism of functors
$\Gamma_\bullet\circ\psi\isom c_G$.
\end{claim}
\begin{pfclaim} By virtue of proposition
\ref{prop_filter-Deligne}(i) we may assume that $(I,\leq)$
is a filtered partially ordered set. Since $I$ is filtered,
the colimit of $\Gamma_\bullet$ commutes with the forgetful
functor $\Phi:\Grp\to\Set$.
Especially, for every $x\in G$ we may find $i(x)\in\Ob(I)$
and $x'\in\Gamma_{i(x)}$ such that $\tau_{i(x)}(x')=x$. Since
$I$ is filtered and $G$ is finite, we may then find
$i\in\Ob(I)$ with $i\geq i(x)$ for every $x\in G$.
Hence the rule : $x\mapsto\Gamma_{i(x),i}(x')$ defines a
set-theoretic section $\sigma_i:G\to\Gamma_i$ of $\tau_i$.
Then for every $j\geq i$, the composition
$\sigma_j:=\Gamma_{ij}\circ\sigma_i$ is a set-theoretic
section of $\tau_j$. After replacing $I$ by the cofinal
subset $\{j\in I~|~j\geq i\}$ we may therefore assume
that there exists a cone
$\sigma_\bullet:c_G\Rightarrow\Phi\circ\Gamma_\bullet$
such that $\tau_i\circ\sigma_i=\one_G$ for every $i\in I$.
Next, for every $x,y\in G$ we may find $i(x,y)\in I$
such that
$$
\sigma_{i(x,y)}(xy)=\sigma_{i(x,y)}(x)\cdot\sigma_{i(x,y)}(y)
$$
(details left to the reader). Again, we may then pick
$i_0\in I$ such that $i_0\geq i(x,y)$ for every $x,y\in G$,
in which case it is easily seen that
$\sigma_{i_0}(xy)=\sigma_{i_0}(x)\cdot\sigma_{i_0}(y)$ for every
$x,y\in G$, {\em i.e.} $\sigma_{i_0}$ is a group homomorphism.
Then clearly $\sigma_j$ is a group homomorphism for every
$j\geq i_0$, so after replacing $I$ by the cofinal subset
of elements $\geq i_0$, we may assume that $\sigma_\bullet$
is even a cone $c_G\Rightarrow\Gamma_\bullet$, and that $I$
admits an initial element $i_0$. Now, let $g_1,\dots,g_k$
be a finite system of generators for $\Gamma_{i_0}$.
For $t=1,\dots,k$ and every $l\in I$ we have
$\tau_l\circ\sigma_l\circ\tau_l(g_t)=\tau_l\circ\Gamma_{i_0,l}(g_t)$;
hence we may find $j_t\in I$ such that
$$
\sigma_{j_t}\circ\tau_{j_t}(g_t)=\Gamma_{i_0,j_t}(g_t).
$$
Then pick again $j\in I$ such that $j\geq j_t$ for every
$t=1,\dots,k$, and set $g'_t:=\Gamma_{i_0,j}(g_t)$ for
every such $t$. It follows that for every $t=1,\dots,k$
we have :
$$
\begin{aligned}
\sigma_j\circ\tau_j(g'_t)&\,=
\sigma_j\circ\tau_j\circ\Gamma_{i_0,j}(g_t) \\
&\,=\sigma_j\circ\tau_{j_t}\circ\Gamma_{i_0,j_t}(g_t) \\
&\,=\Gamma_{j_t,j}\circ\sigma_{j_t}\circ\tau_{j_t}
\circ\Gamma_{i_0,j_t}(g_t) \\
&\,=\Gamma_{j_t,j}\circ\Gamma_{i_0,j_t}(g_t) \\
&\,=g'_t.
\end{aligned}
$$
Notice that $g'_1,\dots,g'_k$ is a generating system
for $\Gamma_j$, since $\Gamma_{i_0,j}$ is a surjective map.
We conclude that $\sigma_j\circ\tau_j=\one_{\Gamma_j}$, so
$\tau_j$ and $\sigma_j$ are mutually inverse group isomorphisms.
After replacing $I$ by the cofinal subset $\{l\in I~|~l\geq j\}$,
we may then assume that $j=i_0$ is the initial element of $I$.
Then, for every $i\in I$ the map
$\sigma_i=\Gamma_{i_0,i}\circ\sigma_{i_0}$ is surjective; but
$\sigma_i$ is also injective by construction, so $\sigma_i$
is an isomorphism for every $i\in I$, and finally, the same
follows for $\tau_i$.
\end{pfclaim}

From claim \ref{cl_fg-groups-are-fp} we deduce that there
exists a morphism $X/h:(X\to uY')\to f$ in $X/u\cC'$ such
that $G_h$ is a finite group, and therefore
$\underline\sigma':=K^nh(\sigma_1,\dots,\sigma_n)\in(K^n)_\sff(Y')$.
Then the class $[\underline\sigma']\in u_!(K^n)(X)$ lies
in $u_!((K^n)_\sff)(X)$ and
$\omega_X([\underline\sigma'])=([\sigma_1],\dots,[\sigma_n])$
as required.
\end{proof}

\begin{theorem}\label{th_CCCP}
In the situation of \eqref{subsec_choo-choo}, suppose that
$u:C\to C'$ is a morphism of sites, and let $K$ be a sheaf
of groups $K$ on $C'$. Then $\tilde u{}^*K$  is a sheaf of
groups on $C$, and we have a natural isomorphism of sheaves :
$$
\omega^{(n)}_{u,K}:\tilde u{}^*((K^n)_\slf)\isom
((\tilde u{}^*K)^n)_\slf
\qquad
\text{for every $n\in\N$}.
$$
\end{theorem}
\begin{proof} By assumption, for every universe $\sV$ such
that $C$ and $C'$ are $\sV$-small, the functor
$\tilde u{}^*:C'^\sim\to C^\sim$ is exact; as in the
proof of proposition \ref{prop_gattiglio}(i), we deduce that
$\tilde u{}^*K$ is a sheaf of groups. We get therefore for
every such $\sV$ a well defined functor
$$
(\tilde u,\sV\tdu\Grp)^*:(C',\Grp)^\sim\to(C,\sV\tdu\Grp)^\sim
\qquad
K\mapsto\tilde u{}^*K.
$$
Next, the unit of adjunction $\eta_K:K\to\tilde u_*\tilde u{}^*K$
induces a morphism of presheaves
$(\eta^n_K)_\sff:(K^n)_\sff\to((\tilde u_*\tilde u{}^*K)^n)_\sff$.
By \eqref{subsec_choo-choo}, we have
$((\tilde u_*\tilde u{}^*K)^n)_\sff=
u^\wedge(((\tilde u{}^*K)^n)_\sff)$, so by adjunction $(\eta^n_K)_\sff$
corresponds to a morphism of presheaves
$u_!((K^n)_\sff)\to((\tilde u{}^*K)^n)_\sff$. After taking associated
sheaves, we get a morphism of sheaves
$\omega:(u_!((K^n)_\sff))^a\to((\tilde u{}^*K)^n)_\slf$. Lastly,
the inclusion of presheaves $(K^n)_\sff\to(K^n)_\slf$ induces an
isomorphism $(u_!((K^n)_\sff))^a\isom\tilde u{}^*((K^n)_\slf)$, by
lemma \ref{lem_cont-funct-site}(ii), and $\omega^{(n)}_{u,K}$ is
the composition of $\omega$ with the inverse of this isomorphism.
With the notation of \eqref{subsec_scampata-bella}, we have
therefore a natural transformation :
$$
\omega^{(n)}_{u,\bullet}:\tilde u{}^*\circ(-)^n_\slf\Rightarrow
(-)^n_\slf\circ(\tilde u,\Grp)^*
\qquad
K\mapsto\omega^{(n)}_{u,K}
$$
and it remains to check that $\omega^{(n)}_{u,K}$ is an isomorphism.
We begin with the following more explicit description of the map
$(\omega^{(n)}_{u,K})_X:\tilde u{}^*((K^n)_\slf)(X)\to
((\tilde u{}^*K)^n)_\slf(X)$, for every $X\in\Ob(\cC)$.
First, by definition every element of $\tilde u{}^*((K^n)_\slf)(X)$
is the equivalence classes $[\underline\sigma]_f$ of a sequence
$\underline\sigma:=(\sigma_1,\dots,\sigma_n)\in K^nY$, for
a morphism $f:X\to uY$ in $\cC$, such that for some covering
sieve $\cS\subset\cC'/Y$, the subgroup
$\La Kg(\sigma_1),\dots,Kg(\sigma_n)\Ra\subset KY'$
is finite for every $(g:Y'\to Y)\in\Ob(\cS)$. Likewise,
every element of $((\tilde u{}^*K)^n)_\slf(X)$ is a sequence
$([\tau_1]_{f_1},\dots,[\tau_n]_{f_n})$ such that
$[\tau_i]_{f_i}\in\tilde u{}^*K(X)$ for every $i=1,\dots,n$
is the equivalence class of a section $\tau_i\in KY_i$, for a
morphism $f_i:X\to uY_i$ in $\cC$, and there exists a covering
sieve $\cT\subset\cC/X$ such that the subgroup
$\La[\tau_1]_{f_1\circ h},\cdots,[\tau_n]_{f_n\circ h}\Ra
\subset(\tilde u{}^*K)(X')$ is finite for every
$(h:X'\to X)\in\Ob(\cT)$.

Now, given such class $[\underline\sigma]_f$ and covering
sieve $\cS$, the sieve $u(\cS)\subset\cC/uY$ generated by
$\{u(g)~|~g\in\Ob(\cS)\}$ covers $uY$ (lemma
\ref{lem_crit-continuity}), and $\cT:=f\times_{uY}u(\cS)$ covers
$X$ for the topology $J$. For every $(h:X'\to X)\in\Ob(\cT)$
there exist $(g:Y'\to Y)\in\Ob(\cS)$ and a morphism
$f':X'\to uY'$ in $\cC$ with $f\circ h=u(g)\circ f'$.
It follows that
$$
[\sigma_i]_{f\circ h}=[\sigma_i]_{u(g)\circ f'}=[Kg(\sigma_i)]_{f'}
\qquad
\text{for $i=1,\dots,n$}
$$
and therefore $\La[\sigma_1]_{f\circ h},\dots,[\sigma_n]_{f\circ h}\Ra$
is a finite subgroup of $(\tilde u{}^*K)(X')$. Then, by unwinding
the definition, we find that $(\omega^{(n)}_{u,K})_X$ is given by
the rule :
$$
[\underline\sigma]_f\mapsto([\sigma_1]_f,\dots,[\sigma_n]_f).
$$
Next, endow $T:=C^\sim$ and $T':=C'^\sim$ with their canonical
topologies, and consider the essentially commutative diagram of
sites provided by lemma \ref{lem_cont-funct-site}(ii) :
$$
{\diagram
T \ar[r]^-{u'} \ar[d]_v \drtwocell\omit{_\beta} &
T' \ar[d]^{v'} \\
C \ar[r]_-u & C'
\enddiagram}\qquad
\text{with\ \ $u':=\tilde u{}^*$,\ \ \ $v':=h^a_{C'}$\ \ and\ \ $v:=h^a_C$}.
$$
After replacing $\sU$ by a larger universe, we may assume
that $T$ and $T'$ are small; we notice :

\begin{claim}\label{cl_ok-for-u-prime}
The theorem holds for $u'$ and every sheaf $H$ on $T'$.
\end{claim}
\begin{pfclaim} Recall that $u'$ is a morphism of lex-sites
(remark \ref{rem_choose-two-univs}(ii)), and in particular
$u'$ fulfills the condition of lemma \ref{lem_slf-and-u_!},
so we get an isomorphism of presheaves
$\omega:u'_!((H^n)_\sff)\isom((u'_!H)^n)_\sff$. Taking into account
lemma \ref{lem_cont-funct-site}(ii), we deduce an isomorphism of
sheaves
$\omega^a:\tilde u{}'^*((H^n)_\slf)\isom((\tilde u{}'^*H)^n)_\slf$.
But a direct inspection of the construction shows that the
latter agrees with $\omega^{(n)}_{u',H}$.
\end{pfclaim}

\begin{claim}\label{cl_ok-also-for-v-and-v-prime}
The theorem holds as well for $v$ and $v'$.
\end{claim}
\begin{pfclaim} It suffices to show the claim for $v$, since
the same argument will work for $v'$ as well. Recall now that
$\breve v{}^*=\tilde v_*:C^\sim\to T^\sim$ is an equivalence
(theorem \ref{th_canon-topos}(iii)), and $\tilde v{}^*$ is
its quasi-inverse. It follows that
$(\tilde v,\Grp)_*:(C,\Grp)^\sim\to(T,\Grp)^\sim$ is also
an equivalence, and its inverse is $(\tilde v,\Grp)^*$.
Then, let $H$ be any sheaf of groups on $C$; in order to
check that $\omega_{v,H}$ is an isomorphism, we may assume
that $H=\tilde v_*L$ for a sheaf of groups $L$ on $T$, and
it suffices to check the commutativity of the diagram :
$$
\xymatrix@C+20pt{ \tilde v{}^*((\tilde v_*L^n)_\slf) \ddouble
\ar[r]^-{\omega_{v,\tilde v_*L}^{(n)}} &
((\tilde v{}^*\tilde v_*L)^n)_\slf \ar[d]^{((\eps_L)^n)_\slf} \\
\tilde v{}^*\tilde v_*((L^n)_\slf) \ar[r]^-{\eps_{(L^n)_\slf}} &
(L^n)_\slf
}$$
where $\eps_\bullet:\tilde v{}^*\tilde v_*\isom\one_{T^\sim}$
is the counit for the adjoint $(\tilde v{}^*,\tilde v_*)$.
However, for every sheaf $F$ on $T$, and every $X\in\Ob(T)$
the elements of $\tilde v{}^*\tilde v_*F(X)$ are the
equivalence classes $[\sigma]_f$, where $f:X\to vY$ is a
morphism in $T$, and $\sigma\in\tilde v_*F(Y)=F(vY)$; with
this notation, the map
$\eps_{F,X}:\tilde v{}^*\tilde v_*F(X)\isom FX$ is given by
the rule : $[\sigma]_f\mapsto Ff(\sigma)$. Then the assertion
follows after simple inspection, by combining with the
foregoing explicit description of $\omega_{v,\tilde v_*L}^{(n)}$.
\end{pfclaim}

We deduce an oriented diagram of categories :
$$
\xymatrix@C+40pt{(T',\Grp)^\sim \ar[rrr]^-{(-)^n_\slf}
\ar[ddd]_{(\tilde u{}',\Grp)^*} & &
\dltwocell\omit{_\ \ \ \ \ \ \ \ \ \ \ \ \ \omega^{(n)}_{v',\bullet}} &
T'^\sim \ar[ddd]^{\tilde u{}'^*} \\
& \dltwocell\omit{^\ \phi} (C',\Grp)^\sim \ar[r]|-{(-)^n_\slf}
\drtwocell\omit{_\ \ \ \ \omega^{(n)}_{u,\bullet}}
\ar[lu]^-{(\tilde v',\Grp)^*} \ar[d]_{(\tilde u,\Grp)^*} &
C'^\sim \ar[d]^{\tilde u{}^*} \ar[ru]_-{\tilde v{}'^*}
\drtwocell\omit{_\ \psi} \\
& (C,\Grp)^\sim \ar[r]|-{(-)^n_\slf}
\drtwocell\omit{_\ \ \ \ \omega^{(n)}_{v,\bullet}}
\ar[ld]_-{(\tilde v,\Grp)^*} & C^\sim \ar[rd]^-{\tilde v^*} & \\
(T,\Grp)^\sim \ar[rrr]_-{(-)^n_\slf} & & & T^\sim.
}$$
which we complete by adding the orientation
$\omega^{(n)}_{u',\bullet}$ for the ``front face''. Here
$\phi$ and $\psi$ denote the identifications induced by
$\beta$ :
$$
(\tilde u{}',\Grp)^*\circ(\tilde v{}',\Grp)^*\isom
(\tilde v,\Grp)^*\circ(\tilde u,\Grp)^*
\qquad\text{and}\qquad
\tilde u{}'^*\circ\tilde v{}'^*\isom\tilde v^*\circ\tilde u^*.
$$
Explicitly, for every sheaf $F'$ on $C'$ and every $X\in\Ob(T)$,
the elements of $\tilde u{}'^*\circ\tilde v{}'^*F(X)$ are the
equivalences classes $[[\tau]_g]_f$, where $f:X\to u'Y$ is a
morphism of $T$, $g:Y\to v'Z$ is a morphism of $T'$, and
$\tau\in FZ$; likewise one describes the elements of
$\tilde v^*\circ\tilde u^*F(X)$, and notice that
$\beta_Z\circ u'(g)\circ f$ is a morphism $X\to vuZ$ in $T$.
Then $\psi_{F,X}:\tilde u{}'^*\circ\tilde v{}'^*F(X)\isom
\tilde v^*\circ\tilde u^*F(X)$ is given by the rule :
$$
[[\tau]_g]_f\mapsto[[\tau]_{\beta_Z\circ u'(g)\circ f}]_{\one_{uZ}}.
$$
If $F$ is a sheaf of groups on $T$, the same rule defines
the group isomorphism $\phi_{F,X}$. Now, a simple inspection
shows that this oriented diagram commutes on $2$-cells, in
the sense of \eqref{subsec_transfer-base-change}. By claims
\ref{cl_ok-for-u-prime} and \ref{cl_ok-also-for-v-and-v-prime}
we know that $\omega^{(n)}_{u',\bullet},\omega^{(n)}_{v,\bullet}$ and
$\omega^{(n)}_{v',\bullet}$ are isomorphisms; we conclude that
the same holds for $\omega^{(n)}_{u,\bullet}$, as stated.
\end{proof}

\begin{corollary}\label{cor_two-assertions}
In the situation of theorem {\em\ref{th_CCCP}}, the following holds :
\begin{enumerate}
\item
If $K$ is ind-finite, the same holds for $\tilde u{}^*K$.
\item
If\/ $\tilde u{}^*$ is a conservative functor and $\tilde u{}^*K$
is ind-finite, then $K$ is ind-finite.
\end{enumerate}
\end{corollary}
\begin{proof} The assertions are immediate consequences of
theorem \ref{th_CCCP} and remark \ref{rem_CCCP}(ii).
\end{proof}

\begin{proposition}\label{prop_ind-finite-and-tildeu_lower-*}
In the situation of \eqref{subsec_choo-choo}, suppose that
$u$ is continuous and there exists a topologically
generating family $G\subset\Ob(\cC')$ such that $uY$ is
quasi-compact for the topology $J$, for every $Y\in G$.
Then, if $H$ is an ind-finite sheaf on $C$, the sheaf\/
$\tilde u_*H$ is ind-finite on $C'$.
\end{proposition}
\begin{proof} Let $X\in\Ob(\cC')$ be any object, and
$\Sigma\subset\tilde u_*H(X)=H(uX)$ any finite set. Pick
a covering family $(f_i:Y_i\to X~|~i\in I)$ with $Y_i\in G$
for every $i\in I$; by remark \ref{rem_ind-finite-and-qc},
the group $\tilde u_*H(Y_i)=H(uY_i)$ is ind-finite, hence
$Hf_i(\Sigma)$ generates a finite subgroup of $H(uY_i)$
for every $i\in I$, whence the assertion.
\end{proof}

\sset\subsubsection{}\label{subsec_automorph-cart-sect}
$C:=(\cC,J)$ be a site, and $\cA\to\cC$ a fibration;
for every $X\in\Ob(\cC)$ and every cartesian section
$\sigma\in\Ob(\cA(X))$ we define the presheaf on $\cC/X$ :
$$
\cCart^\times(\sigma):=\cCart(\sigma^\times,\sigma^\times)
$$
(notation of \eqref{subsec_up-to-gpds-in-fibs}). Clearly,
for every $f\in\Ob(\cC/X)$, the set $\cCart^\times(\sigma)(f)$
carries a natural group structure, given by composition of
automorphisms of the cartesian section
$(\sigma\circ f_*)^\times$ of $\cA^\times$, so
$\cCart^\times(\sigma)$ is a presheaf of groups on $\cC/X$.

\begin{definition} With the notation of
\eqref{subsec_automorph-cart-sect}, we say that $\cA$ is
an {\em ind-finite prestack} on $C$, if $\cCart^\times(\sigma)$
an ind-finite sheaf of groups on $C/X$, for every $X\in\Ob(\cC)$,
and every $\sigma\in\Ob(\cA(X))$ (notation of
\eqref{sec_Localization-topoi}). We say that $\cA$ is an
{\em ind-finite stack} on $C$, if it is  both an ind-finite
prestack and a stack on $C$.
\end{definition}

\sset\subsubsection{}\label{subsec_Van-Gogh}
Let $u:C:=(\cC,J)\to C':=(\cC',J')$ be a morphism of small
sites, $\cE$ a stack on $C'$; let also $X\in\Ob(\cC')$,
and $\sigma,\tau\in\Ob(\cE(X))$. Set $\cE':=\St(u)^*\cE$,
and denote
$$
\eta_\cE:\cE\to\St(u)_*\cE'
$$
the unit of the $2$-adjunction for the $2$-adjoint
pair $(\St(u)^*,\St(u)_*)$. According to
\eqref{subsec_pull-back-presh-of-cart}, the natural
projection $\pi:\St(u)_*\cE'\to\cE'$ induces an
equivalence of categories :
$$
u^*_{|X}:\cE'(uX)\isom\St(u)_*\cE'(X).
$$
We may then find $\sigma',\tau'\in\Ob(\cE'(uX))$ and
isomorphisms
$$
\omega_1:\eta_\cE\circ\sigma\isom u^*_{|X}(\sigma')
\qquad
\omega_2:\eta_\cE\circ\tau\isom u^*_{|X}(\tau').
$$
Taking into account lemma \ref{lem_crit-0-1_separation}
we deduce isomorphisms of sheaves on $C'/X$ :
$$
\cCart(\eta_\cE\circ\sigma,\eta_\cE\circ\tau)\isom
\cCart(u^*_{|X}\sigma',u^*_{|X}\tau')\isom
\tilde u_{|X*}\cCart(\sigma',\tau')
$$
namely the composition of the isomorphism \eqref{eq_pariah}
induced by $\omega_1$ and $\omega_2$, with the isomorphism
of presheaves given in \eqref{subsec_pull-back-presh-of-cart}.
Then \eqref{eq_F-pariah} yields as well a morphism of sheaves :
$$
\cCart(\sigma,\tau)\to
\cCart(\eta_\cE\circ\sigma,\eta_\cE\circ\tau)
\qquad
\beta\mapsto\eta_\cE*\beta.
$$
By adjunction, the composition of these morphisms induces
a morphism of sheaves
on $C/uX$ :
\set\begin{equation}\label{eq_bravo}
\tilde u{}^*_{|X}\cCart(\sigma,\tau)\to\cCart(\sigma',\tau').
\end{equation}

\begin{remark}
(i)\ \
Let $(f:Y\to uX)\in\Ob(\cC/uX)$. The evaluation at $f$
of the morphism \eqref{eq_bravo} can be described explicitly
as follows. For every element
$[\alpha]\in\tilde u{}^*_{|X}\cCart(\sigma,\tau)(f)$ we may find
a covering family $(k_\lambda:Y_\lambda\to Y~|~\lambda\in\Lambda)$
and for every $\lambda\in\Lambda$ a morphism
$h_\lambda:Y_\lambda\to uZ_\lambda$ in $\cC$ and a morphism
$g_\lambda:Z_\lambda\to X$ in $\cC'$ with
$f\circ k_\lambda=u(g_\lambda)\circ h_\lambda$, such that $[\alpha]$
is represented by a system of morphisms of cartesian sections
$(\alpha_\lambda:\sigma\circ g_{\lambda*}\Rightarrow\tau\circ g_{\lambda*})$.
Set
$$
\alpha^*_\lambda:=\pi(\omega_{2,g_\lambda}\circ
\eta_\cE(\alpha_{\lambda,\one_{Z_\lambda}})\circ\omega_{1,g_\lambda}^{-1}):
\sigma'(ug_\lambda)\to\tau'(ug_\lambda)
\qquad
\text{for every $\lambda\in\Lambda$}.
$$
Then for every $\lambda\in\Lambda$ there exists a unique morphism
$\alpha^{**}_\lambda:\sigma'(f\circ k_\lambda)\to\tau'(f\circ k_\lambda)$
in $\cE'_{Y_\lambda}$ that makes commute the following diagram in $\cE'$ :
$$
\xymatrix{ \sigma'(f\circ k_\lambda) \ar[r]^-{\alpha^{**}_\lambda}
\ar[d]_{\sigma'(h_\lambda/uX)} &
\tau'(f\circ k_\lambda) \ar[d]^{\tau'(h_\lambda/uX)} \\
\sigma'(ug_\lambda) \ar[r]^-{\alpha^*_\lambda} & \tau'(ug_\lambda)
}$$
and there exists a unique morphism of cartesian sections
$\alpha^{**}:\sigma'\circ f_*\Rightarrow\tau'\circ f*$ such
that $\alpha^{**}_{k_\lambda}=\alpha^{**}_\lambda$ for every
$\lambda\in\Lambda$. Then the map \eqref{eq_bravo} is
given by the rule : $[\alpha]\mapsto\alpha^{**}$.

(ii)\ \
The morphism \eqref{eq_bravo} depends obviously on the
choice of $\sigma'$, $\tau'$, $\omega_1$ and $\omega_2$.
However, say that $\sigma'',\tau''\in\Ob(\cE'(uX))$ are
two other cartesian sections with isomorphisms
$\omega'_1:\eta_\cE\circ\sigma\isom u^*_{|X}(\sigma'')$
and $\omega'_2:\eta_\cE\circ\tau\isom u^*_{|X}(\tau'')$.
Then there exists unique isomorphisms
$\rho_1:\sigma'\isom\sigma''$ and $\rho_2:\tau'\isom\tau''$
in $\cE'(uX)$ such that
$$
\omega'_i=u^*_{|X}(\rho_i)\circ\omega_i
\qquad
\text{for $i=1,2$}
$$
and in light of (i) it is easily seen that we get a
commutative diagram :
$$
\xymatrix{ & \cCart(\sigma,\tau) \ar[ld] \ar[rd] \\
\cCart(\sigma',\tau') \ar[rr]^-\sim & & \cCart(\sigma'',\tau'') 
}$$
whose downward arrows are the morphisms \eqref{eq_bravo}
relative to the choices $(\omega_1,\omega_2)$ and respectively
$(\omega'_1,\omega'_2)$, and whose horizontal arrow is the
isomorphism \eqref{eq_pariah} induced by $(\rho_1,\rho_2)$.
\end{remark}

\begin{proposition}\label{prop_long-proof}
The morphism \eqref{eq_bravo} is an isomorphism.
\end{proposition}
\begin{proof} Let us check that \eqref{eq_bravo} is a
monomorphism. Thus, let $(f:Y\to uX)\in\Ob(\cC/uX)$, and
$[\alpha],[\beta]\in\tilde u{}^*_{|X}\cCart(\sigma,\tau)(f)$
whose images agree in $\cCart(\sigma',\tau')(f)$; we need
to show that $[\alpha]=[\beta]$. We may find a covering
family $(k_\lambda:Y_\lambda\to Y~|~\lambda\in\Lambda)$
such that the restrictions $[\alpha_\lambda],[\beta_\lambda]
\in\tilde u{}^*_{|X}\cCart(\sigma,\tau)(f\circ k_\lambda)$
of $[\alpha]$ and $[\beta]$ are in the image of the natural
map $u_{|X!}\cCart(\sigma,\tau)(f\circ k_\lambda)\to
\tilde u{}^*_{|X}\cCart(\sigma,\tau)(f\circ k_\lambda)$, for every
$\lambda\in\Lambda$. Then the images of $[\alpha_\lambda]$ and
$[\beta_\lambda]$ agree in $\cCart(\sigma',\tau')(f\circ k_\lambda)$,
and it suffices to check that $[\alpha_\lambda]=[\beta_\lambda]$
for every such $\lambda$. Thus, we may assume from start that
$[\alpha],[\beta]$ are in the image of the natural map
$u_{|X!}\cCart(\sigma,\tau)(f)\to
\tilde u{}^*_{|X}\cCart(\sigma,\tau)(f)$. This means that there
exists a commutative diagram
$$
\xymatrix{ uZ' \ar[rd]_{u(g')} &
Y \ar[d]^f \ar[r]^-h \ar[l]_-{h'} &
uZ \ar[ld]^{u(g)} \\
& uX
}$$
in $\cC$ such that $[\alpha]$ and $[\beta]$ are represented by
natural $\cC$-transformations
$$
\alpha:\sigma\circ g_*\Rightarrow\tau\circ g_*
\qquad
\beta:\sigma\circ g'_*\Rightarrow\tau\circ g'_*.
$$
By proposition \ref{prop_morph-of-sites}, the fibration
$\ss:\cC/u\cC'\to\cC$ is locally cofiltered relative to
the topology $J$ of $C$, hence there exist a covering
family $(k_\lambda:Y_\lambda\to Y~|~\lambda\in\Lambda)$ and for
every $\lambda\in\Lambda$ a morphism
$h''_\lambda:Y_\lambda\to uZ_\lambda$ in $\cC$, and two morphisms
$l_\lambda:Z_\lambda\to Z$, $l'_\lambda:Z_\lambda\to Z'$ such that
$$
h\circ k_\lambda=u(l_\lambda)\circ h''_\lambda
\qquad\text{and}\qquad
h'\circ k_\lambda=u(l'_\lambda)\circ h''_\lambda
\qquad
\text{for every $\lambda\in\Lambda$}.
$$
It suffices to check that the images of $[\beta]$ and
$[\alpha]$ in $\cCart(\sigma,\tau)(f\circ k_\lambda)$
coincide for every $\lambda\in\Lambda$; but these images
are the classes of $[\alpha*l_\lambda]$ and respectively
$[\beta*l'_\lambda]$. Hence, we may assume from start that
$Z=Z'$, $h=h'$ and $g=g'$. With the notation of
\eqref{subsec_pull-back-presh-of-cart}, the condition
on $[\alpha]$ and $[\beta]$ amounts then to the identity :
\set\begin{equation}\label{eq_woppity}
u^*_{|Z}(\eta_\cE*\alpha)*h_*=u^*_{|Z}(\eta_\cE*\beta)*h_*
\qquad\text{in $\cE'(Y)$}.
\end{equation}
Recall that the evaluation functor $\sev^\cE_Z:\cE(Z)\to\cE_Z$
is an equivalence of categories for every $Z\in\Ob(\cC)$,
natural with respect to morphism in $\cC$ (claim
\ref{cl_split-fibrations}); this equivalence assigns to the
cartesian sections $\sigma\circ g_*$ and $\tau\circ g_*$
their evaluations $\sigma_g:=\sigma(g)$ and
$\tau_g:=\tau(g)$ in $\sc_X$, and to $\alpha$ and $\beta$
the morphisms $\alpha_0:=\alpha_{\one_Z}:\sigma_g\to\tau_g$
and $\beta_0:=\beta_{\one_Z}:\sigma_g\to\tau_g$ in $\cE_Z$.
We wish next to similarly interpret the
identity \eqref{eq_woppity} as an equality of morphisms
in the fibre category $\cE'_Y$. To this aim, let us pick
a unital cleavage for $\cE$, and let $\sc$ be its associated
unital pseudo-functor; denote also by $\gamma^\sc$ the
coherence constraint of $\sc$. Recall that for every
$Y\in\Ob(\cC)$, the fibre category $\Fib(u)_!\cE_Y$ of
the $\cC$-fibration $\Fib(u)_!\cE$ represents the
$2$-colimit of the pseudo-functor
$$
F_Y:=\sc\circ\st_Y^o:(Y/u\cC')^o\to\bCat
\qquad
(g:Y\to uZ)\mapsto\sc_Z=\cE_Z.
$$
The latter is realized by the localization
$\cFib(F_Y)[\Sigma_Y^{-1}]$, where $\Sigma_Y$ denotes the set of
cartesian morphisms of $\cFib(F_Y)$. Recall that the objects
of $\cFib(F_Y)[\Sigma_Y^{-1}]$ are the data
$$
(Z,h:Y\to uZ,T)
\qquad
\text{with $(Z,h)\in\Ob(Y/u\cC')$ and $T\in\Ob(\cE_Z)$}.
$$
The morphisms are the pairs $[Y/g,t]:(Z,h,T)\to(Z',h',T')$
where $Y/g:h\to h'$ is a morphism in $Y/u\cC'$ and $t:T\to\sc_gT'$
is a morphism in $\cE_Z$. Moreover, every morphism
$k:Y\to Y'$ in $\cC$ induces a $\cC$-cartesian functor
$$
k^*:\cFib(F_{Y'})\to\cFib(F_Y)
\qquad
(Z,h,T)\mapsto(Z,h\circ k,T)
\qquad
[Y'/g,t]\mapsto[Y/g,t]
$$
which therefore extends uniquely to a functor
$$
\cFib(F_{Y'})[\Sigma^{-1}_{Y'}]\to\cFib(F_Y)[\Sigma^{-1}_Y]
$$
and the system of such functors provides a natural split
cleavage for $\Fib(u)_!\cE$. Notice then the rules :
$Y\mapsto\cFib(F_Y)$ and $k\mapsto k^*$ for every
$Y\in\Ob(\cC)$ and every morphism $k$ of $\cC$, define
a pseudo-functor $\sd:\cC^o\to\bCat$, and the foregoing
description amounts to saying that $\Fib(u)_!\cE$ is
represented by the fibration which in the proof of
proposition \ref{prop_go-for-kill} is denoted
$$
\cFib(\sd)[\Sigma^{-1}_\bullet]\to\cC.
$$
By proposition \ref{prop_go-for-kill}, the natural
functor $j:\Fib(u)_!\cE\to\cE'$ factors then as a
composition
$$
\Fib(u)_!\cE\xrightarrow{L_{\{\Sigma\}}}
\cFib(\sd)\{\Sigma^{-1}_\bullet\}\xrightarrow{j'}\cE'.
$$
Furthermore, for every $Z\in\Ob(\cC')$ we have a natural
identification $u^*_{|uZ}:\St(u)_*\cE'_Z\isom\cE'_{uZ}$. Hence,
$\eta_\cE$ restricts to a functor of fibre categories
$\eta_{\cE,Z}:\cE_Z\to\cE'_{uZ}$ for every such $Z$.
However, $\cE_Z$ is also the fibre category over
$\one_{uZ}\in\Ob(uZ/u\cC')$ of the $(uZ/u\cC')$-fibration
$\cFib(F_{uZ})$; in terms of the foregoing realization of
$\Fib(u)_!\cE_{uZ}$, the functor $\eta_{\cE,Z}$ is then just
given by the rule :
$$
\cE_Z\to\cFib(F_{uZ})[\Sigma^{-1}_{uZ}]\xrightarrow{j_{uZ}}\cE'_{uZ}
\qquad
T\mapsto j_{uZ}(Z,\one_{uZ},T)
$$
where $j_{uZ}$ is the restriction of $j$. Summing up,
we find that the objects $\eta_\cE(\sigma_g)$ and
$\eta_\cE(\tau_g)$ of $\cE'_{uZ}$ are respectively
$j_{uZ}(Z,\one_{uZ},\sigma_g)$ and $j_{uZ}(Z,\one_{uZ},\tau_g)$,
and \eqref{eq_woppity} translates as the identity:
$$
j_{uZ}[Y/\one_Z,\alpha_0]=j_{uZ}[Y/\one_Z,\beta_0].
$$
But since $j$ is $2$-covering (proposition
\ref{prop_unit-sep-cov-faith}(i)), we deduce that there
exists a covering family
$(k_\lambda:Y_\lambda\to Y~|~\lambda\in\Lambda)$ in $\cC$ such
that $[Y_\lambda/\one_Z,\alpha_0]=[Y_\lambda/\one_Z,\beta_0]$
in $\cFib(F_{Y_\lambda})[\Sigma^{-1}_{Y_\lambda}]$ for every
$\lambda\in\Lambda$. Since it suffices to check that the
images of $[\alpha]$ and $[\beta]$ coincide in
$\cCart(\sigma,\tau)(f\circ k_\lambda)$ for every $\lambda$,
we may then assume from start that
$[Y/\one_Z,\alpha_0]=[Y/\one_Z,\beta_0]$ in
$\cFib(F_Y)[\Sigma^{-1}_Y]$.

\begin{claim}\label{cl_yes-we-win}
The system $\Sigma_\bullet:=(\Sigma_Y~|~Y\in\Ob(\cC))$ admits
a right local calculus of fractions.
\end{claim}
\begin{pfclaim} Conditions (LCF1) and (LCF2) of definition
\ref{def_LCF} obviously hold for $\Sigma_\bullet$. Next,
let $(Y/k',l):(h',T')\to(h,T)$ be a morphism in $\cFib(F_Y)$
and $(Y/k'',t):(h'',T'')\to(h,T)$ an element of $\Sigma_Y$;
{\em i.e.} we have a commutative diagram of morphisms of $\cC$ :
$$
\xymatrix{ uZ' \ar[rd]_{u(k')} &
Y \ar[d]^h \ar[r]^-{h''} \ar[l]_-{h'} &
uZ'' \ar[ld]^{u(k'')} \\
& uZ
}$$
together with a morphism $l:T'\to\sc_{k'}T$ and an isomorphism
$t:T''\isom\sc_{k''}T$ in $\cA_Y$. By proposition
\ref{prop_morph-of-sites}, the fibration $\ss:\cC/u\cC'\to\cC$
is locally cofiltered relative to the topology $J$ of $C$, hence
there exist a covering family
$(g_\lambda:Y_\lambda\to Y~|~\lambda\in\Lambda)$ and for
every $\lambda\in\Lambda$ a morphism
$h_\lambda:Y_\lambda\to uZ_\lambda$ in $\cC$, and two morphisms
$k'_\lambda:Z_\lambda\to Z'$, $k''_\lambda:Z_\lambda\to Z''$ such that
$$
h'\circ g_\lambda=u(k'_\lambda)\circ h_\lambda
\qquad
h''\circ g_\lambda=u(k''_\lambda)\circ h_\lambda
\qquad
k'\circ k'_\lambda=k''\circ k''_\lambda
\qquad
\text{for every $\lambda\in\Lambda$}
$$
whence a well defined morphism in $\cA_{Y_\lambda}$, for every
$\lambda\in\Lambda$ :
$$
l_\lambda:\sc_{k'_\lambda}T'\xrightarrow{\sc_{k'_\lambda}l}
\sc_{k'_\lambda}\sc_{k'}T\xrightarrow{\gamma^\sc_{(k'_\lambda,k'),T}}
\sc_{k'\circ k'_\lambda}T\xrightarrow{\gamma^{\sc\,-1}_{(k''_\lambda,k''),T}}
\sc_{k''_\lambda}\sc_{k''}T\xrightarrow{(\sc_{k''_\lambda}t)^{-1}}
\sc_{k''_\lambda}T''
$$
which then yield a commutative diagram in $\cFib(F_{Y_\lambda})$
for every such $\lambda$ :
$$
\xymatrix@C+40pt{ (h_\lambda,\sc_{k'_\lambda}T')
\ar[r]^-{(Y_\lambda/k'_\lambda,t_\lambda)}
\ar[d]_{(Y_\lambda/k''_\lambda,l_\lambda)} &
g_\lambda^*(h',T') \ar[d]^{g_\lambda^*(Y/k',l)} \\
g_\lambda^*(h'',T'') \ar[r]^-{g_\lambda^*(Y/k'',t)} &
g_\lambda^*(Y/h,T)
}$$
where we take for $t_\lambda$ the identity of $\sc_{k'_\lambda}T'$.
This shows that condition (LCF3) holds for $\Sigma_\bullet$.
Lastly, suppose that we have a commutative diagram in
$\cFib(\cA_Y)$ :
$$
\xymatrix@C+30pt{ (h,T) \ar@<.5ex>[r]^-{(Y/k,l)}
\ar@<-.5ex>[r]_-{(Y/k',l')} & (h',T') \ar[r]^-{(Y/k'',t)} &
(h'',T'')
}$$
with $(Y/k'',t)\in\Sigma_Y$. This means that we have
a commutative diagram in $\cC$ :
$$
\xymatrix@C+20pt{
& Y \ar@/_1pc/[ld]_h \ar[d]^{h'} \ar@/^1pc/[rd]^{h''} \\
uZ \ar@<.5ex>[r]^-{u(k)} \ar@<-.5ex>[r]_-{u(k')} &
uZ' \ar[r]^-{u(k'')} & uZ''
}$$
and morphisms $l:T\to\sc_kT'$, $l':T\to\sc_{k'}T'$ and
an isomorphism $t:T'\isom\sc_{k''}T''$ in $\cA_Y$ with
$$
k''\circ k=k''\circ k'
\qquad\text{and}\qquad
\gamma^\sc_{(k,k''),T''}\circ\sc_k(t)\circ l=
\gamma^\sc_{(k',k''),T''}\circ\sc_{k'}(t)\circ l'.
$$
Since the fibration $\ss:\cC/u\cC'\to\cC$ is locally cofiltered,
there exist a covering family
$(g_\lambda:Y_\lambda\to Y~|~\lambda\in\Lambda)$ and for every
$\lambda\in\Lambda$, morphisms $h_\lambda:Y_\lambda\to uZ_\lambda$
in $\cC$ and $k_\lambda:Z_\lambda\to Z$ in $\cC'$ such that
$h\circ g_\lambda=u(k_\lambda)\circ h_\lambda$ in $\cC$
and $k\circ k_\lambda=k'\circ k_\lambda$ in $\cC'$. We compute :
$$
\begin{aligned}
\gamma^\sc_{(k\circ k_\lambda,k''),T''}\circ\sc_{k\circ k_\lambda}(t)
\circ\gamma^\sc_{(k,k_\lambda),T'}\circ\sc_{k_\lambda}(l)
&=\gamma^\sc_{(k\circ k_\lambda,k''),T''}\circ
\gamma^\sc_{(k,k_\lambda),\sc_{k''}T''}
\circ\sc_{k_\lambda}\sc_k(t)\circ\sc_{k_\lambda}(l) \\
&=\gamma^\sc_{(k_\lambda,k''\circ k),T''}\circ
\sc_{k_\lambda}(\gamma^\sc_{(k,k''),T''})
\circ\sc_{k_\lambda}\sc_k(t)\circ\sc_{k_\lambda}(l) \\
&=\gamma^\sc_{(k_\lambda,k''\circ k'),T''}\circ
\sc_{k_\lambda}(\gamma^\sc_{(k,k''),T''}\circ\sc_k(t)\circ l) \\
&=\gamma^\sc_{(k_\lambda,k''\circ k'),T''}\circ
\sc_{k_\lambda}(\gamma^\sc_{(k',k''),T''}\circ\sc_{k'}(t)\circ l') \\
&=\gamma^\sc_{(k_\lambda,k''\circ k'),T''}\!\circ\!
\sc_{k_\lambda}(\gamma^\sc_{(k',k''),T''})\!\circ\!
\sc_{k_\lambda}\sc_{k'}(t)\!\circ\!\sc_{k_\lambda}(l') \\
&=\gamma^\sc_{(k'\circ k_\lambda,k''),T''}\circ
\gamma^\sc_{(k',k_\lambda),\sc_{k''}T''}\circ
\sc_{k_\lambda}\sc_{k'}(t)\circ\sc_{k_\lambda}(l') \\
&=\gamma^\sc_{(k'\circ k_\lambda,k''),T''}\circ\sc_{k'\circ k_\lambda}(t)
\circ\gamma^\sc_{(k',k_\lambda),T'}\circ\sc_{k_\lambda}(l')
\end{aligned}
$$
whence $\gamma^\sc_{(k,k_\lambda),T'}\circ\sc_{k_\lambda}(l)=
\gamma^\sc_{(k',k_\lambda),T'}\circ\sc_{k_\lambda}(l')$. We deduce
a commutative diagram in $\cFib(F_{Y_\lambda})$
$$
\xymatrix@C+40pt{
(h_\lambda,\sc_{k_\lambda}T) \ar[r]^-{(Y/k_\lambda,t_\lambda)} &
g_\lambda^*(h,T) \ar@<.5ex>[r]^-{g^*_\lambda(Y/k,l)}
\ar@<-.5ex>[r]_-{g^*_\lambda(Y/k',l')} &
g_\lambda^*(h',T')
}\qquad\text{for every $\lambda\in\Lambda$}
$$
where $t_\lambda$ is the identity of $\sc_{k_\lambda}T$. This
shows that (LFC4) holds for $\Sigma_\bullet$.
\end{pfclaim}

To ease notation, set $A:=L_{\{\Sigma\}}(Z,h,\sigma_g)$
and $B:=L_{\{\Sigma\}}(Z,h,\tau_g)$, and define the
presheaf $M_{A,B}$ on $\cC/uZ$ and the sheaf $H_{A,B}$ on
$C/uZ$ as in \eqref{subsec_here-we-go-at-last}. Then
$$
L_{\{\Sigma\}}[Y/\one_Z,\alpha]=
L_{\{\Sigma\}}[Y/\one_Z,\beta]\in H_{A,B}(h)
$$
so there exists a covering family
$(k_\lambda:Y_\lambda\to Y~|~\lambda\in\Lambda)$ such that
the images of $[Y_\lambda/\one_Z,\alpha]$ and
$[Y_\lambda/\one_Z,\beta]$ agree in $M_{A,B}(h\circ g_\lambda)$.
Hence as usual we may assume from start that the images
of $[Y/\one_Z,\alpha_0]$ and $[Y/\one_Z,\beta_0]$ agree
in $M_{A,B}(h)$. Then, in view of claim \ref{cl_yes-we-win},
condition (b) of \eqref{subsec_here-we-go-at-last} tells
us that there exist a covering family
$(k_\lambda:Y_\lambda\to Y~|~\lambda\in\Lambda)$, and for
every $\lambda\in\Lambda$ a morphism
$h_\lambda:Y_\lambda\to uZ_\lambda$ in $\cC$ and a morphism
$l_\lambda:Z_\lambda\to Z$ in $\cC'$ such that
$h\circ k_\lambda=u(l_\lambda)\circ h_\lambda$ and with
$\sc_{l_\lambda}(\alpha_0)=\sc_{l_\lambda}(\beta_0)$. This
means that the images of $[\alpha]$ and $[\beta]$ agree
in $\tilde u{}^*_{|X}\cCart(\sigma,\tau)(f\circ k_\lambda)$
for every $\lambda\in\Lambda$, so $[\alpha]=[\beta]$, as
required.

Next, let us check that \eqref{eq_bravo} is an epimorphism.
To this aim, set $\cE'':=\Fib(u)_!\cE$, and let
$\eta'_\cE:\cE\to\Fib(u)^*\cE''$ be the unit of the
$2$-adjunction for the $2$-adjoint pair $(\Fib(u)_!,\Fib(u)^*)$;
as in\eqref{subsec_Van-Gogh}, the projection
$\pi':\Fib(u)^*\cE''\to\cE''$ induces an equivalence
of categories
$$
u^*_{|X}:\cE''(uX)\isom\Fib(u)^*\cE''(X)
$$
so we may find $\sigma'',\tau''\in\Ob(\cE''(uX))$ with
isomorphisms :
$$
\eta'_\cE\circ\sigma\isom u^*_{|X}(\sigma'')
\qquad
\eta'_\cE\circ\tau\isom u^*_{|X}(\tau'').
$$
Since $\Fib(u)^*(j)\circ\eta'_\cE:\cE\to\St(u)_*\cE'$ is
isomorphic to $\eta_\cE$, we deduce isomorphisms
$$
\sigma'\isom j\circ\sigma''
\qquad
\tau'\isom j\circ\tau''
$$
and it suffices to show that the corresponding morphism
of sheaves on $C/uX$
$$
\tilde u{}^*_{|X}\cCart(\sigma,\tau)\to
\cCart(j\circ\sigma'',j\circ\tau'')
$$
is an epimorphism. Thus, let $(f:Y\to X)\in\Ob(\cC/X)$ and
$\mu'\in\cCart(j\circ\sigma'',j\circ\tau'')(f)$; we need to
show that
$\mu':j\circ\sigma''\circ f_*\Rightarrow j\circ\tau''\circ f_*$ is
the image of some $[\mu]\in\tilde u{}^*_{|X}\cCart(\sigma,\tau)(f)$.
Set $\sigma''_0:=\sigma''(f)$, $\tau''_0:=\tau''(f)$, and
recall that $\mu'$ is determined by its evaluation
$$
\mu'_0:=\mu'_{\one_Y}:j(\sigma''_0)\to j(\tau''_0).
$$
Since $j$ is $1$-covering (proposition
\ref{prop_unit-sep-cov-faith}(i)), there exist a covering
family $(k_\lambda:Y_\lambda\to Y~|~\lambda\in\Lambda)$, and
for every $\lambda\in\Lambda$ a morphism
$\mu''_\lambda:\sigma''(f\circ k_\lambda)\to
\tau''(f\circ k_\lambda)$ in $\cE''_{Y_\lambda}$ such that
$\mu'_{k_\lambda}=j(\mu''_\lambda)$. Now, suppose that
for every $\lambda\in\Lambda$ there exists $\mu_\lambda
\in\tilde u{}^*_{|X}\cCart(\sigma,\tau)(f\circ k_\lambda)$
whose image in $\cCart(\sigma',\tau')(f\circ k_\lambda)$
equals $\mu'*k_\lambda$. Then, for every
$\lambda,\lambda'\in\Lambda$, and every pair of morphisms
$Y_{\lambda'}\xleftarrow{k''_{\lambda\lambda'}}Y_{\lambda\lambda'}
\xrightarrow{k'_{\lambda\lambda'}}Y_\lambda$ in $\cC$ such
that $k_{\lambda\lambda'}:=
k_\lambda\circ k'_{\lambda\lambda'}=k_{\lambda'}\circ k''_{\lambda\lambda'}$,
let $\mu_{\lambda\lambda'}$ and $\mu_{\lambda'\lambda}$ be the
images of $\mu_\lambda$ and respectively $\mu_{\lambda'}$
in $\tilde u^*_{|X}\cCart(\sigma,\tau)(k_{\lambda\lambda'}\circ f)$;
it follows that $\mu_{\lambda\lambda'}$ and $\mu_{\lambda'\lambda}$
are both mapped to $\mu'*k_{\lambda\lambda'}$ by the morphism
\eqref{eq_bravo}. Since we already know that the latter
is a monomorphism, we deduce that
$\mu_{\lambda\lambda'}=\mu_{\lambda'\lambda}$ for every such pair
of morphisms $k'_{\lambda\lambda'}, k''_{\lambda\lambda'}$. Then the
system $(\mu_\lambda~|~\lambda\in\Lambda)$ determines
a unique $\mu\in\tilde u^*_{|X}\cCart(\sigma,\tau)(f)$
whose image in $\cCart(\sigma',\tau')(f)$ equals $\mu'$,
as sought. Thus, we may replace $\mu$ by $\mu*k_\lambda$
for every $\lambda\in\Lambda$, and assume from start that
there exists a morphism
$$
\mu''_0:L_{\{\Sigma\}}\sigma''_0\to L_{\{\Sigma\}}\tau''_0
\qquad
\text{in $\cFib(\sd)\{\Sigma^{-1}_\bullet\}_Y$ such that}
\qquad
\mu'_0=j'(\mu''_0).
$$
Set $\sigma''':=L_{\{\Sigma\}}\sigma''$ and
$\tau''':=L_{\{\Sigma\}}\tau''$; then $\mu''_0$ corresponds
to a unique morphism of cartesian sections
$\mu'':\sigma'''\circ f_*\Rightarrow\tau'''\circ f_*$
with $j'*\mu''=\mu'$. Recall that $\sigma''_0$ (resp.
$\tau''_0$) is a datum $(Z,h,T)$ (resp. $(Z',h',T')$).
Now, according to condition (a) of
\eqref{subsec_here-we-go-at-last}, we may find
a covering family
$(k_\lambda:Y_\lambda\to Y~|~\lambda\in\Lambda)$, and for
every $\lambda\in\Lambda$, an object
$(Z_\lambda,h_\lambda:Y_\lambda\to uZ_\lambda,T_\lambda)$ of
$\cFib(F_{Y_\lambda})$, and morphisms $(Z,h\circ k_\lambda,T)
\xleftarrow{(Y_\lambda/l_\lambda,s_\lambda)}(Z_\lambda,h_\lambda,T_\lambda)
\xrightarrow{(Y_\lambda/l'_\lambda,t_\lambda)}(Z',h'\circ k_\lambda,T')$
in $\cFib(F_{Y_\lambda})$, where $s_\lambda:T_\lambda\to\sc_{l_\lambda}T$
is an isomorphism in $\cE_{Z_\lambda}$, and such that we have a
commutative diagram in $\cFib(\sd)\{\Sigma^{-1}_\bullet\}_{Y_\lambda}$ :
$$
\xymatrix@C+120pt{ (Z,h\circ k_\lambda,T) \ar[d]
\ar[r]^-{L_{\{\Sigma\}}(Y_\lambda/l'_\lambda,t_\lambda)\circ
L_{\{\Sigma\}}(Y_\lambda/l_\lambda,s_\lambda)^{-1}} &
(Z',h'\circ k_\lambda,T) \ar[d] \\
\sigma'''(f\circ k_\lambda)
\ar[r]^-{\mu''_{k_\lambda}} & \tau'''(f\circ k_\lambda)
}$$
whose vertical arrows are isomorphisms. The latter
is equivalent to a commutative diagram :
$$
\xymatrix@C+80pt{ (Z_\lambda,h_\lambda,\sc_{l_\lambda}T) \ar[d]
\ar[r]^-{L_{\{\Sigma\}}(Y_\lambda/\one_{Z_\lambda},t_\lambda\circ s^{-1}_\lambda)}
& (Z_\lambda,h_\lambda,\sc_{l'_\lambda}T') \ar[d] \\
\sigma'''(f\circ k_\lambda)
\ar[r]^-{\mu''_{k_\lambda}} &
\tau'''(f\circ k_\lambda)
}$$
whose vertical arrows are again isomorphisms. After
replacing $\sigma'''$ and $\tau'''$ by isomorphic
cartesian sections, we may therefore assume that
$\mu''_{k_\lambda}=
L_{\{\Sigma\}}(Y_\lambda/\one_{Z_\lambda},t_\lambda\circ s^{-1}_\lambda)$
for every $\lambda\in\Lambda$.

Arguing as in the foregoing, we may then assume that
there exist a morphism $h:Y\to uZ$ in $\cC$ and a morphism
$t:T\to T'$ in $\cE_Z$ such that
$\mu'_0=j(Y/\one_Z,t):j(Z,h,T)\to j(Z,h,T')$.

Next, since the fibration $\ss:\cC/u\cC'\to\cC$ is locally
cofiltered over the site $C$, we may find a covering family
$(k_\lambda:Y_\lambda\to Y~|~\lambda\in\Lambda)$, and for every
$\lambda\in\Lambda$ a morphism $h_\lambda:Y_\lambda\to uZ_\lambda$
in $\cC$ and two morphisms
$X\xleftarrow{g_\lambda}Z_\lambda\xrightarrow{g'_\lambda}Z$
such that $f\circ k_\lambda=g_\lambda\circ h_\lambda$
and $h''\circ k_\lambda=g'_\lambda\circ h_\lambda$. Then
again, we may replace $\sigma'$ and $\tau'$ by isomorphic
cartesian sections, and assume that
$\mu'_{k_\lambda}=j(Y_\lambda/\one_{Z_\lambda},\sc_{g'_\lambda}t):
(Z_\lambda,h_\lambda,\sc_{g'_\lambda}T)\to
(Z_\lambda,h_\lambda,\sc_{g'_\lambda}T')$. Arguing once again
as in the foregoing, we may finally assume that we have
still $\mu'_0=j(Y/\one_Z,t):j(Z,h,T)\to j(Z,h,T')$,
and moreover there exists a morphism $g:Z\to X$ in $\cC$
such that $u(g)\circ h=f$. But then, $t$ corresponds to
a morphism of cartesian sections
$\mu:\sigma\circ f_*\to\tau\circ f_*$, and the pair
$(Y/g:f\to h,\mu)$ yields the sought section $[\mu]$.
\end{proof}

\begin{theorem}\label{th_ind-fin-and-pullback-stacks}
In the situation of \eqref{subsec_Van-Gogh}, the following holds :
\begin{enumerate}
\item
If $\cE$ is ind-finite, then the same holds for $\cE'$.
\item
If $\tilde u^*$ is a conservative functor and $\cE'$ is
ind-finite, then $\cE$ is ind-finite.
\end{enumerate}
\end{theorem}
\begin{proof}(i): In view of proposition \ref{prop_Van-Gogh},
we may replace $\cE$ by $\cE^\times$, and assume from start
that $\cE$ is a stack in groupoids, in which case we need
to show that for every $Y\in\Ob(\cC)$ and every
$\sigma'\in\Ob(\cE'(Y))$, the sheaf of groups
$\cCart(\sigma',\sigma')$ is ind-finite on $C/Y$.

Since the natural functor $j:\cE'':=\Fib(u)_!\cE\to\cE'$ is
$0$-covering (proposition \ref{prop_unit-sep-cov-faith}(i)),
there exist a covering family
$(\phi_\lambda:Y_\lambda\to Y~|~\lambda\in\Lambda)$, and
for every $\lambda\in\Lambda$ a cartesian section
$\sigma''_\lambda\in\Ob(\cE''(Y_\lambda))$ with an isomorphism
$j\circ\sigma''_\lambda\isom\sigma'\circ\phi_{\lambda*}$.
Clearly it suffices to show that
$\cCart(\sigma'\circ\phi_{\lambda*},\sigma'\circ\phi_{\lambda*})$
is ind-finite for every $\lambda\in\Lambda$. We may thus
replace $Y$ by $Y_\lambda$ for every such $\lambda$, and
assume from start that $\sigma'=j\circ\sigma''$ for some
$\sigma''\in\Ob(\cE''(Y))$. Recall that
$\sigma''(\one_Y)\in\Ob(\cE''_Y)$ is a datum $(X,f:Y\to uX,T)$
where $(X,f)\in\Ob(Y/u\cC')$ and $T\in\Ob(\cE_X)$.
Let then $\tau:\cC'/X\to\cE$ and $\tau':\cC/uX\to\cE''$
be cartesian sections with $\tau_{\one_X}=T$ and
$\tau'(\one_{uX}):=(X,\one_{uX},T)$. As explained in the
proof of of proposition \ref{prop_long-proof}, we have
an isomorphism :
$$
\eta'_\cE\circ\tau\isom u^*_{|X}(\tau')
$$
where $\eta'_\cE:\cE\to\Fib(u)^*\cE''$ is the unit of a
$2$-adjunction for the $2$-adjoint pair $(\Fib(u)_!,\Fib(u)^*)$.
Since $\Fib(u)^*(j)\circ\eta'_\cE:\cE\to\St(u)_*\cE'$ is
isomorphic to the unit $\eta_\cE:\cE\to\St(u)_*\cE'$ of a
$2$-adjunction for the $2$-adjoint pair $(\Fib(u)_!,\Fib(u)^*)$,
there follows an isomorphism
$$
\eta_\cE\circ\tau\isom u^*_{|X}(j\circ\tau').
$$
From proposition \ref{prop_long-proof} and corollary
\ref{cor_two-assertions}(i) we conclude that
$\cCart(j\circ\tau',j\circ\tau')$ is ind-finite on $C/uX$.
On the other hand, we have $\tau'\circ f_*\simeq\sigma''$,
whence $j\circ\tau'\circ f_*\simeq\sigma'$. Combining with
\eqref{eq_triviality}, there follows an isomorphism of
sheaves $\cCart(\sigma',\sigma')\isom
(f_*)^\sim_*\cCart(j\circ\tau',j\circ\tau')$. But recall
that $f_*$ is both continuous and cocontinuous for the
topologies of the sites $C/Y$ and $C/uX$ (remark
\ref{rem_continue-local}(i)); then the assertion follows
from proposition \ref{prop_gattiglio}(ii) and lemma
\ref{lem_improve}(i).

(ii): We may reduce as in (i) to the case where $\cE$
is a stack in groupoids. Next notice that if $\tilde u^*$
is conservative, $\tilde u^*_{|X}$ is conservative for
every $X\in\Ob(\cC')$, by virtue of proposition
\ref{prop_localize-continuity}(ii). Then the assertion
follows immediately from proposition \ref{prop_long-proof}
and corollary \ref{cor_two-assertions}(ii).
\end{proof}

\begin{proposition}\label{prop_ind-finite-and-lower-*-stacks}
In the situation of \eqref{subsec_choo-choo}, suppose that
$u$ is a weak morphism of sites, and moreover there exists
a topologically generating family $G\subset\Ob(\cC')$ such
that $uY$ is quasi-compact for the topology $J$, for every
$Y\in G$. Then, for every ind-finite stack $\cE$ on $C$, the
stack $\St(u)_*\cE$ is ind-finite on $C'$.
\end{proposition}
\begin{proof} In view of the pseudo-commutative diagram
\eqref{eq_burkini} we may assume that $\cE$ is a stack in
groupoids, in which case the same holds for $\cE':=\St(u)_*\cE$.
Now, let $X\in\Ob(\cC')$ and $\sigma,\tau\in\cE'(X)$. Recall
that $u$ induces a continuous functor $u_{|X}:\cC'/X\to\cC/uX$
(proposition \ref{prop_localize-continuity}(i)), and moreover
$G/X\subset\Ob(\cC'/X)$ is a topologically generating family,
according to \eqref{subsec_topol-on-C-over-X}; then it is easily
seen that for every $(Y\xrightarrow{g}X)\in\Ob(G/X)$ the object
$uY\xrightarrow{u(g)}uX$ of $\cC/uX$ is quasi-compact.
On the other hand, arguing as in the proof of proposition
\ref{prop_C0-C4}(i), we find $\sigma',\tau'\in\cE(uX)$ and
isomorphisms of sheaves
$$
\cCart(\sigma,\tau)\isom\tilde u_{|X*}\cCart(\sigma',\tau').
$$
The assertion then follows from proposition
\ref{prop_ind-finite-and-tildeu_lower-*}.
\end{proof}

\sset\subsubsection{Stalk of a stack over a $T$-point}
\label{subsec_stalks-of-stacks}
Let $T$ be a topos, and $\xi:=(\xi^*,\xi_*):\Set\to T$ a $T$-point
(see definition \ref{def_T-point}(i)). For every stack $\cE$ on
the canonical site $\Can(T)$ we set
$$
\cE_\xi:=\St(\xi^*)^*\cE.
$$
Likewise, for every morphism $\phi:\cE\to\cF$ of stacks on
$\Can(T)$, we let
$$
\phi_\xi:=\St(\xi^*)^*(\phi):\cE_\xi\to\cF_\xi.
$$

\begin{remark}
(i)\ \
Recall that $\Set$ is equivalent to $(\bone,J)^\sim$,
where $\bone$ is the category with one object and one
morphism, and $J$ is the unique topology on $\bone$.
It follows easily that a stack on the canonical site
$\Can(\Set)$ of $\Set$ amounts to the datum of a category,
and a morphism of stacks on $\Can(\Set)$ is just an
arbitrary functor.

(ii)\ \
Likewise, a functor $F:\cA\to\cB$ between categories is
$0$-covering (resp. $1$-covering, resp. $2$-covering)
when regarded as a morphism of stacks on $\Can(\Set)$,
if and only if it is essentially surjective (resp. full,
resp. faithful) : the details shall be left to the reader.
\end{remark}

\begin{proposition}\label{prop_stalks-for-stacks}
With the notation of \eqref{subsec_stalks-of-stacks}, let
$\Omega$ be a conservative set of $T$-points (see definition
{\em\ref{def_T-point}(iv)}). Then for every $i=0,1,2$ the
following conditions are equivalent :
\begin{enumerate}
\alphaenu
\item
$\phi$ is $i$-covering.
\item
The morphism $\phi_\xi$ is $i$-covering for every $\xi\in\Omega$.
\end{enumerate}
\end{proposition}
\begin{proof}(a)$\Rightarrow$(b) follows from proposition
\ref{prop_Fib-lower-shriek-and-u}(ii).

(b)$\Rightarrow$(a): For $i=1$ and $i=2$, the assertion
follows easily from lemma \ref{lem_charact-1-2-coverings}
and proposition \ref{prop_long-proof}. Suppose then that
$\phi_\xi$ is $0$-covering for every $\xi\in\Omega$; by
lemma \ref{lem_parthenon}(i), it follows that
$\pi_0(\phi_\xi):\pi_0(\cE_\xi)\to\pi_0(\cF_\xi)$ is a
surjection. But according to \eqref{subsec_pull-back-pi_0^C},
$\pi_0(\cE_\xi)$ is naturally identified with the stalk
$\pi_0^T(\cE)_\xi$ of the sheaf $\pi_0^T(\cE)$ on $\Can(T)$,
and likewise for $\pi_0(\cF_\xi)$. Thus, the morphism of
sheaves $\pi_0^T(\phi):\pi_0^T(\cE)\to\pi_0^T(\cF)$ is an
epimorphism, and therefore $\phi$ is $0$-covering, again
by lemma \ref{lem_parthenon}(i).
\end{proof}

\subsection{Stacks on fibred sites and fibred topoi}
\label{sec_C2-for-fibred-sites}
Let $\cC$ and $\cC'$ be two categories, $u:\cC'\to\cC$ any
functor, $(\cA,p,J_\bullet)$ a fibred site over $\cC$, and
$(\cA',p'_\bullet,J'_\bullet):=\cC'\times_\cC(\cA,p,J_\bullet)$
the induced fibred site as in \eqref{subsec_pullback-fib-site};
denote also by $(\cA,J_\cA)$ and $(\cA',J_{\cA'})$ the respective
fibred sites. Then it follows easily from proposition
\ref{prop_sheaves-on-tot-site}(ii,iv) that the projection
$\pi:\cA'\to\cA$ fulfills condition (C2) of
\eqref{subsec_weak-strong-continuity}, relative to the
topologies $J_\cA$ and $J_{\cA'}$, so it is a weak morphism
of sites $\pi:(\cA,J_\cA)\to(\cA',J_{\cA'})$, and we get a
well defined pseudo-functor
$$
\St(\pi)_*:\Stack(\cA,J_\cA)\to\Stack(\cA',J_{\cA'}).
$$
If $\cC$ and $\cC'$ are small and $(\cA_X,J_X)$ is a $\sU$-site
for every $X\in\Ob(\cC)$, then $(\cA,J_\cA)$ and $(\cA',J_{\cA'})$
are also $\sU$-sites (remark \ref{rem_fibred-site}(iii)), and
in this case also the left $2$-adjoint $\St(\pi)^*$ of
$\St(\pi)_*$ is well defined; moreover, $\pi$ is cocontinuous
(see \eqref{subsec_pullback-fib-site}), so we have a
pseudo-natural equivalence $\St(\pi)_*\isom\St(\breve\pi)^*$
(proposition \ref{prop_breve-for-stacks}(ii)). 

As a special case, for every $X\in\Ob(\cC)$ the inclusion
functor $i_X:\cA_X\to\cA$ is a weak morphism of sites
$(\cA,J_\cA)\to(\cA_X,J_X)$, and we get a well defined
pseudo-functor
$$
\St(i_X)_*:\Stack(\cA,J_\cA)\to\Stack(\cA_X,J_X).
$$
If $\cC$ is small and $(\cA_X,J_X)$ is a $\sU$-site for every
$X\in\Ob(\cC)$, also $\St(i_X)^*$ is well defined.

\begin{proposition}\label{prop_fibrewise-i-cov}
In the situation of \eqref{sec_C2-for-fibred-sites}, we have :
\begin{enumerate}
\item
Let $\phi:\cE\to\cE'$ be a cartesian functor of fibrations over
$\cA$, and $i\in\{0,1,2\}$. The following conditions are
equivalent :
\begin{enumerate}
\item
$\phi$ is $i$-covering for the topology $J_\cA$.
\item
$\Fib(i_X)^*(\phi)$ is $i$-covering for the topology
$J_X$, for every $X\in\Ob(\cC)$.
\end{enumerate}
\item
Suppose that $\cC$ is small, and $(\cA_X,J_X)$ is a $\sU$-site
for every $X\in\Ob(\cC)$. Then the following diagram of
$2$-categories is pseudo-commutative :
$$
{\spreaddiagramrows{-5pt}\diagram
\Fib(\cA) \ar[rr]^-{\Fib(i_X)^*} \ar[d]_{(-)^a} & &
\Fib(\cA_X) \ar[d]^{(-)^a} \\
\Stack(\cA,J_\cA) \ar[rr]^-{\St(i_X)_*} & & \Stack(\cA_X,J_X)
\enddiagram}
\qquad
\text{for every $X\in\Ob(\cC)$}.
$$
\end{enumerate}
\end{proposition}
\begin{proof}(i.a)$\Rightarrow$(i.b): This follows directly
from propositions \ref{prop_cocont-and-stacks}(i) and
\ref{prop_sheaves-on-tot-site}(iii).

(i.b)$\Rightarrow$(i.a): Suppose first that $i=0$, and let
$X\in\Ob(\cA)$ and $E'\in\Ob(\cE'_X)$; pick cleavages
for $\cE$ and $\cE'$, and let $\sc$ and $\sc'$ be the
corresponding associated pseudo-functors. By assumption,
there exist a covering family
$f_\bullet:(f_i:X_i\to X~|~i\in I)$ for the topology $J_{pX}$,
and for every $i\in I$ an object $E_i\in\Ob(\cE_{X_i})$ such
that $\phi E_i$ is isomorphic to $\sc'_{f_i}E'$ in
$\cE'_{X_i}$. Since $f_\bullet$ is also a covering family
for the topology $J_\cA$, this shows that $\phi$ is
$0$-covering, as required. The cases where $i=1,2$ are
similar : the details shall be left to the reader.

(ii): The assertion is a special case of corollary
\ref{cor_weak-morph-of-sites-and-cocont}.
\end{proof}

\begin{corollary}\label{cor_lims-and-totsites}
Let $j\leq 2$ be any integer. The following holds :
\begin{enumerate}
\item
In the situation of \eqref{subsec_cofilt-lims-lexsites}, let
$\cE$ be a fibration on $\cC$. Then $\cE$ is a $j$-separated
prestack on $C$ if and only if\/ $\Fib(\pi_i)^*(\cE)$ is a
$j$-separated prestack on $C_i$ for every $i\in\Ob(I)$.
\item
In the situation of \eqref{sec_C2-for-fibred-sites}, let
$\cE$ be a fibration on $\cA$. Then $\cE$ is a $j$-separated
prestack on $(\cA,J_\cA)$ if and only if\/ $\Fib(i_X)^*(\cE)$ is
$j$-separated on $(\cA_X,J_X)$ for every $X\in\Ob(\cC)$.
\end{enumerate}
\end{corollary}
\begin{proof}(i): The condition is necessary, by virtue of
corollary \ref{cor_another-condition} and proposition
\ref{prop_C0-C4}(i,ii). For the converse, we argue as in
the proof of proposition \ref{prop_2-lims-lex-sites}(ii).
Indeed, suppose that $\Fib(\pi_i)^*(\cE)$ is a $j$-separated
prestack on $C_i$ for every $i\in\Ob(I)$; according to remark
\ref{rem_sheaves}(ii) there exists a topology $J'$ on $\cC$
such that, for every $X\in\Ob(\cC)$, the set $J'(X)$ consists
of the sieves $\cS\subset\cC/X$ of universal $\cE$-$j$-descent.
It then suffices to prove that $J\subset J'$. We come down
to checking the following. For every $i\in\Ob(I)$, every
$X\in\Ob(\cC_i)$, every covering family
$g_\bullet:=(g_\lambda:X_\lambda\to X~|~\lambda\in\Lambda)$
for the topology $J_i$, and every morphism $f:Y\to\pi_iX$
in $\cC$, the family $(Y\times_{\pi_iX}\pi_i(g_\lambda):
Y\times_{\pi_iX}\pi_iX_\lambda\to Y~|~\lambda\in\Lambda)$ generates
a sieve of $\cE$-$j$-descent. Arguing as in the proof of proposition
\ref{prop_2-lims-lex-sites}(ii), we are first reduced to the case
where $Y=\pi_iY'$ for some $Y'\in\Ob(\cC_i)$ and $f=\pi_i(t)$ for
some morphism $t:Y'\to X$ of $\cC_i$, and then we may even replace
the family $g_\bullet$ by $(Y'\times_Xg_\lambda~|~\lambda\in\Lambda)$
and reduce further to the case where $t=\one_X$. Finally, pick a
cleavage $\sc$ for $\cE$; we need to show that the natural functor
$$
\cE(\pi_iX)\to\Desc(\cE,\pi_i(g_\bullet),\sc)
$$
is $j$-faithful (notation of \eqref{subsec_descnt-data}). But
since $\pi_i$ is left exact, this functor is naturally identified
with the corresponding functor for $\cE_i:=\Fib(\pi_i)^*(\cE)$ :
$$
\cE_i(X)\to\Desc(\cE_i,g_\bullet,\sc\circ\pi_i^o)
$$
and the latter is $j$-faithful, since $\cE_i$ is a $j$-separated
prestack on $C_i$.

(ii): Again, the condition is necessary, by corollary
\ref{cor_another-condition} and proposition
\ref{prop_C0-C4}(i,ii). For the converse we argue as in the
proof of proposition \ref{prop_sheaves-on-tot-site}.
Namely, we let $J'$ be the topology on $\cA$ such that for
every $A\in\Ob(\cA)$ the set $J'(A)$ consists of the sieves
$\cS\subset\cA/A$ of universal $\cE$-$j$-descent, and we check
that $J_\cA\subset J'$, assuming that $\cE_X:=\Fib(i_X)^*(\cE)$ is a
$j$-separated prestack on $(\cA_X,J_X)$ for every $X\in\Ob(\cC)$.
To this aim, we consider any morphism $f:A'\to A$ in $\cA$ and
any family $(g_\lambda:B_\lambda\to A~|~\lambda\in\Lambda)$ of
morphisms in $\cA_{pA}$ covering $A$ for the topology $J_{pA}$;
after choosing a cleavage for $\cA$, we then construct the family
$g'_\bullet:=(g'_\lambda:B'_\lambda\to A'~|~\lambda\in\Lambda)$ of
morphisms of $\cA_{pA'}$ as in the proof of proposition
\ref{prop_sheaves-on-tot-site}. Let $i_{pA'}:\cA_{pA'}\to\cA$
be the inclusion functor; taking into account lemma
\ref{lem_was-pair-of-sieves}(ii), we reduce to checking that
the sieve of $\cA/A'$ generated by the family
$i_{pA'}(g'_\bullet):=(i_{pA'}(g_\lambda)~|~\lambda\in\Lambda)$
is of $\cE$-$j$-descent. Now, choose a cleavage $\sc$ for
$\cE$, and set
$$
B'_{\lambda\mu}:=B'_\lambda\times_{A'}B'_\mu
\qquad\text{and}\qquad
B'_{\lambda\mu\nu}:=B'_{\lambda\mu}\times_{A'}B'_\nu
\qquad\text{for every $\lambda,\mu,\nu\in\Lambda$}
$$
and notice that $i_{pA'}B'_{\lambda\mu}$ represents the fibre
product $i_{pA'}B'_\lambda\times_{i_{pA'}A'}B'_\mu$, and likewise
$i_{pA'}B'_{\lambda\mu\nu}$ represents
$i_{pA'}B'_{\lambda\mu}\times_{i_{pA'}A'}i_{pA'}B'_\nu$ (remark
\ref{rem_fibred-site}(ii)), so that the category
$\Desc(\cE,i_{pA'}(g'_\bullet),\sc)$ is well defined, and
we are reduced to checking that the induced map
$$
\cE(i_{pA'}A')\to\Desc(\cE,i_{pA'}(g'_\bullet),\sc)
$$
is $j$-faithful. But the latter is naturally identified
with the analogous map
$$
\cE_{pA'}(A')\to\Desc(\cE_{pA'},g'_\bullet,\sc\circ i^o_{pA'})
$$
which is $j$-faithful by assumption.
\end{proof}

\begin{proposition}\label{prop_no-hand-waving}
In the situation of \eqref{subsec_box-a-louer}, suppose that
$(\cA_{i,X},J_{i,X})$ is a $\sU$-site for every $X\in\Ob(\cC)$
and $i=0,1$. Then the induced base change transformation :
$$
\Upsilon(\St(\one_{u\circ\pi_1})_*^\gamma):
\St(v)^*\circ\St(\pi_1)_*\to\St(\pi_0)_*\circ\St(u)^*
$$
is a pseudo-natural equivalence.
\end{proposition}
\begin{proof} Recall that the projections $\pi_0$ and $\pi_1$
satisfy condition (C2) from \eqref{subsec_weak-strong-continuity} :
see \eqref{sec_C2-for-fibred-sites}; in particular, they are
weak morphisms of sites, and moreover they are cocontinuous
(see \eqref{subsec_pullback-fib-site}). We are then in the
situation contemplated in \eqref{subsec_rotate-the-cube}, and
by corollary \ref{cor_conditions-for-trivial-bc}(ii), both
$\Delta^{\pi_0}$ and $\Delta^{\pi_1}$ are pseudo-natural equivalences,
and the same holds for $\Upsilon(\Delta^u)$ and $\Upsilon(\Delta^v)$,
by virtue of corollary \ref{cor_conditions-for-trivial-bc}(iii).
Hence, arguing as in remark \ref{rem_transit-base-change}(i),
we are reduced to checking that
$\Upsilon((\one^\sim_{u\circ\pi_1},\bCat)_*)$ is a pseudo-natural
equivalence. The latter assertion follows from proposition
\ref{prop_Devos} and the discussion of
\eqref{subsec_from-beta-dagger-to-tilde}.
\end{proof}

\begin{definition}\label{def_weak-fibred-morph}
Let $(\cA,p,J_\bullet)$ and $(\cA',p',J'_\bullet)$ be two fibred
sites over a category $\cC$. A {\em weak morphism of fibred
sites} $\phi:(\cA,p,J_\bullet)\to(\cA',p',J'_\bullet)$ is a
functor $\phi:\cA'\to\cA$ with $p\circ\phi=p'$, whose
restriction $\phi_X:\cA'_X\to\cA_X$ is a weak morphism of
sites $(\cA_X,J_X)\to(\cA'_X,J'_X)$, for every $X\in\Ob(\cC)$.
\end{definition}

\begin{proposition}\label{prop_fibred-weak-is-weak}
Let $(\cA,p,J_\bullet)$ and $(\cA',p',J'_\bullet)$ be two fibred
lex-sites over a category $\cC$, and
$\phi:(\cA,p,J_\bullet)\to(\cA',p',J'_\bullet)$ a weak morphism
of fibred sites. Then $\phi$ is a weak morphism of the respective
total sites $(\cA,J_\cA)\to(\cA',J_{\cA'})$.
\end{proposition}
\begin{proof} For every $X\in\Ob(\cC)$, let $i_x:\cA_X\to\cA$
and $i'_X:\cA'_X\to\cA'$ be the inclusion functors, and
$\phi_X:\cA'_X\to\cA_X$ the restriction of $\phi$; in light
of corollary \ref{cor_lims-and-totsites}(ii), it suffices to
check that for every stack $\cE$ over $(\cA,J_\cA)$, the fibration
$\cF:=\Fib(i'_X)^*\circ\Fib(\phi)^*(\cE)$ is a stack on
$(\cA'_X,J'_X)$. But we have a natural identification
$\cF\isom\Fib(\phi_X)^*\circ\Fib(i_X)^*(\cE)$, and $\Fib(i_X)^*(\cE)$
is a stack on $(\cA_X,J_X)$, again by corollary
\ref{cor_lims-and-totsites}(ii), whence the contention.
\end{proof}

\sset\subsubsection{}\label{subsec_localization-morph-of-sites}
Let $(\cA,p,J_\bullet)$ be a small fibred lex-site over a small
category $\cC$; pick a unital cleavage $\blambda$ for $p$, and
let $\sc:\cC^o\to\bCat$ be the associated pseudo-functor. Clearly
$\sc^o$ factors through a pseudo-functor
$$
\tilde\sc:\cC\to\lex.\Site
\qquad
X\mapsto(\cA_X,J_X)
$$
and the forgetful functor $\lex.\Site\to\bCat^o$.
Let $\Sigma_\cA$ be the set of cartesian morphisms of $\cA$.
By remark \ref{ex_filter-2-colim-in-Cat}(i), we know that
$\cA[\Sigma_\cA^{-1}]$ represents the strong $2$-colimit of $\sc$ :
$$
\Pscolim{\cC^o}\sc\isom\cA[\Sigma_\cA^{-1}].
$$
{\em Suppose now that $\cC$ is cofiltered}. Then, the
construction of \eqref{subsec_cofilt-lims-lexsites} endows
$\cA[\Sigma^{-1}_\cA]$ with a natural topology $J^*_\cA$ such that
$(\cA[\Sigma^{-1}_\cA],J^*_\cA)$ is a strong $2$-limit of $\tilde\sc$,
and again by direct inspection we see that the localization
functor $L_\cA:\cA\to\cA[\Sigma_\cA^{-1}]$ is continuous for the
topology $J^*_\cA$ and the topology $J_\cA$ of the total site
$(\cA,J_\cA)$. Moreover, for every $X\in\Ob(X)$, the composition
$\cA_X\to\cA[\Sigma^{-1}_\cA]$ of the inclusion functor
$i_X:\cA_X\to\cA$ and $F_\cA$, is a morphism of sites
$$
L_\cA\circ i_X:(\cA[\Sigma^{-1}_\cA],J^*_\cA)\to(\cA_X,J_X)
$$
and the rule : $X\mapsto L_\cA\circ i_X$ defines both a universal
pseudo-cone
$\tilde\tau:\sF_{(\cA[\Sigma^{-1}_\cA],J^*_\cA)}\Rightarrow\tilde\sc$
and a universal pseudo-cocone
$\tau:\sc\Rightarrow\sF_{\cA[\Sigma^{-1}_\cA]}$. Recall as well that
$\Sigma_\cA$ admits a right calculus of fraction (example
\ref{ex_filter-2-colim-in-Cat}(ii)), so that the category
$X/L_\cA\cA$ is cofiltered for every $X\in\Ob(\cA[\Sigma_\cA^{-1}])$
(proposition \ref{prop_rcf-and-cofilt-commas}(i)). Then the source
fibration $\ss:\cA[\Sigma^{-1}_\cA]/L_\cA\cA\to\cA[\Sigma^{-1}_\cA]$
is trivially locally cofiltered, and from proposition
\ref{prop_morph-of-sites} we see that $L_\cA$ is a morphism of sites
$$
L_\cA:(\cA[\Sigma^{-1}_\cA],J^*_\cA)\to(\cA,J_\cA).
$$

\begin{proposition}\label{prop_local-morph-of-sites}
In the situation of \eqref{subsec_localization-morph-of-sites},
the following holds :
\begin{enumerate}
\item
The functor
$L^\sim_{\cA*}:(\cA[\Sigma_\cA^{-1}],J^*_\cA)^\sim\to(\cA,J_\cA)^\sim$
is fully faithful.
\item
Suppose that every object of $\cA$ is quasi-compact for the
topology $J_\cA$ (see definition {\em\ref{def_qcoh-obj-site}}),
and let $F$ and $\cE$ be respectively a sheaf and a stack
on $(\cA,J_\cA)$. Then:
\begin{enumerate}
\item
$L_{\cA!}F$ is a sheaf on $(\cA[\Sigma^{-1}_\cA],J^*_\cA)$.
\item
$\Fib(L_\cA)_!\cE$ is a stack on $(\cA[\Sigma^{-1}_\cA],J^*_\cA)$.
\end{enumerate}
\end{enumerate}
\end{proposition}
\begin{proof}(i) follows easily from corollary \ref{cor_of-localization}.

Next, let $F\in\Ob(\cA^\wedge)$ and $A\in\Ob(\cA[\Sigma^{-1}_\cA])$;
recall that $L_{\cA!}F(A)$ represents the colimit of
$$
F\circ\st^o_A:(A/L_\cA\cA)^o\to\Set
\qquad
(A\to L_\cA B)\mapsto FB.
$$
We define as follows a functor $\Phi_A:\cC/pA\to A/L_\cA\cA$.
For every $\phi\in\Ob(\cC/pA)$ we let
$$
\Phi_A(\phi):=(\sc_\phi A,(L_\cA\blambda(A,\phi))^{-1}:
A\to L_\cA(\sc_\phi A))
$$
(see \eqref{subsec_choice-of-pseudo}). For every morphism
$\psi/pA:\phi\to\phi'$ of $\cC/pA$ we have $\phi=\phi'\circ\psi$,
and we set
$$
\Phi_A(\psi/pA):=A/L_\cA(\lambda(\one_A,\psi)):
\Phi_A(\phi)\to\Phi_A(\phi').
$$

\begin{claim}\label{cl_Phi-final}
The functor $\Phi_A$ is final.
\end{claim}
\begin{pfclaim} Since $\cC$ is cofiltered, the same holds
for $\cC/A$ (example \ref{ex_filtered-final}(i)), hence we
may apply the criterion of (the dual of) lemma
\ref{lem_filtered-final}(i). First, for every object
$(B,f:A\to L_\cA B)$ of $A/L_\cA\cA$ we may find a morphism
$\phi:pB\to pA$ of $\cC$ and a morphism
$f':\sc_\phi A\to B$ in $p^{-1}(pB)$ such that
$f=L_\cA(f')\circ L_\cA(\lambda(A,\phi))^{-1}$ (details left
to the reader), so that $f'$ yields a morphism
$\Phi_A(\phi)\to(B,f)$ in $A/L_\cA\cA$. Next, let $\phi:X\to pA$
be any object of $\cC/pA$ and
$A/g,A/h:\Phi_A(\phi)\to(B,f:A\to L_\cA B)$ two morphisms of
$A/L_\cA\cA$; by definition, we have
$g=f\circ L_\cA(\lambda(A,\phi))=h$, so condition (b) of lemma
\ref{lem_filtered-final}(i) is trivially verified.
\end{pfclaim}

By claim \ref{cl_Phi-final}, the set $L_{\cA!}F(A)$ also
represents the colimit of the functor
$$
F\circ(\st_A\circ\Phi_A)^o:(\cC/pA)^o\to\Set
\qquad
(\phi:X\to pA)\mapsto F(\sc_\phi A).
$$

(ii.a): Notice first that, under the stated assumptions, for
every $X\in\Ob(\cC)$, every $B\in\Ob(\cA_X)$ is quasi-compact
for the topology $J_X$ : indeed, if $\cS\in J_X(B)$, let $\cS'$
be the sieve of $\cA/B$ generated by $\cS$; we have
$\cS'\in J_\cA(B)$, so by assumption, there exists a finite
subset $\{f_i:B_i\to B~|~i\in I\}\subset\cS'$ that generates
a sieve $\cS''$ covering $B$ in $(\cA,J_\cA)$. But for every
$i\in I$ there exist $(f'_i:C_i\to B)\in\cS$ and a morphism
$g_i:B_i\to C_i$ in $\cA$ such that $f_i=f'_i\circ g_i$; the
family $(f'_i~|~i\in I)$ still generates a covering sieve
for $J_\cA$, and then the same family covers $B$ also in the
site $(\cA_X,J_X)$. Next we remark :

\begin{claim}\label{cl_cofinal-sieves}
Let $A\in\Ob(\cA[\Sigma^{-1}_\cA])=\Ob(\cA)$, and $\cS$ a
sieve of $\cA[\Sigma^{-1}_\cA]/A$. Then $\cS\in J^*_\cA(A)$
if and only if there exist a cartesian morphism
$h:A'\to A$ in $\cA$ and a finite covering family
$(f_i:A'_i\to A')$ in $(\cA_{pA'},J_{pA'})$
with $L_\cA(h\circ f_i)\in\cS$ for every $i\in I$.
\end{claim}
\begin{pfclaim} Since $L_\cA$ is continuous and $L_\cA h$ is
an isomorphism in $\cA[\Sigma^{-1}_\cA]$, the stated condition
implies that $\cS\in J^*_\cA(A)$. For the converse, recall
that $J^*_\cA$ is the coarsest topology such that $L_\cA\circ i_X$
is a continuous functor for the topology $J_X$ and $J^*_\cA$,
for every $X\in\Ob(\cC)$. Thus, for every $B\in\Ob(\cA)$,
denote by $J^{**}_\cA(B)$ the set of all sieves of
$\cA[\Sigma^{-1}_\cA]/B$ fulfilling the condition of the claim;
it suffices then to check that
$J^{**}_\cA:=(J^{**}_\cA(B)~|~B\in\Ob(A))$ is a topology on
$\cA[\Sigma^{-1}_\cA]$, and that $L_\cA\circ i_X$ is a
continuous functor for the topologies $J_X$ and $J^{**}_\cA$,
for every $X\in\Ob(\cC)$. To check stability under base change
for $J^{**}_\cA$ (see definition \ref{def_topology}(i)), consider
any morphism $g:B\to A$ in $\cA[\Sigma^{-1}_\cA]$; we need to
show that $\cS\times_Ag\in J^{**}_\cA(B)$. However, since
$\Sigma_\cA$ admits a right calculus of fraction, there
exist a cartesian morphism $h':B'\to B$ in $\cA$ and a
morphism $k:B'\to A$ in $\cA$ such that
$g=L_\cA(k)\circ L_\cA(h')^{-1}$. Suppose now that
$\cS':=\cS\times_AL_\cA(k)\in J^{**}_{B'}$; this means that
there exist a cartesian morphism $h'':B''\to B'$ and a finite
covering family $(f'_i:B''_i\to B')$ in $(\cA_{pB'},J_{pB'})$
such that $L_\cA(h''\circ f'_i)\in\cS'$ for every
$i\in I$. Therefore
$L_\cA(h'\circ h''\circ f'_i)\in\cS'\times_{B'}L_\cA(h')^{-1}=
\cS\times_Ag$, and since $h'\circ h''$ is cartesian, we
get $\cS\times_Ag\in J^{**}_\cA$. Thus, we may replace $B$
by $B'$, and assume from start that $g=L_\cA(k)$ for some
morphism $k:B\to A$ in $\cA$. In this case, since $\cC$
is cofiltered, we can find $X\in\Ob(\cC)$ with morphisms
$\phi:X\to pB$ and $\psi:X\to pA'$ such that
$p(g)\circ\phi=p(h)\circ\psi$; then we may also find a
morphism $k':B'\to A''$ in $p^{-1}X$, and cartesian morphisms
$h':B'\to B$ and $h'':A''\to A'$ in $\cA$ such that $p(h')=\phi$,
$p(h'')=\psi$, and $k\circ h'=h\circ h''\circ k'$ in $\cA$. For
every $i\in I$ there exist a morphism $f'_i:A''_i\to A''$ in
$p^{-1}X$ and a cartesian morphism $h''_i:A''_i\to A'_i$
such that $f_i\circ h''_i=h''\circ f'_i$. Since $\sc_\psi$
is continuous for the topologies $J_{pA'}$ and $J_X$, it
follows easily that $(f'_i~|~i\in I)$ is a covering
family in $(\cA_X,J_X)$ (details left to the reader),
and by construction $L_\cA(h\circ h''\circ f'_i)\in\cS$
for every $i\in I$. Since $h\circ h''$ is cartesian, we
may therefore replace $A'$ by $A''$, $h$ by $h\circ h''$
and the family $(f_i~|~i\in I)$ by $(f'_i~|~i\in I)$, and
assume from start that there exist a morphism $k':B'\to A'$
in $\cA_{pA'}$ and a cartesian morphism $h':B'\to B$ in
$\cA$ such that $h\circ k'=k\circ h'$. Let
$(f'_i:A'_i\times_{A'}B'\to B'~|~i\in I)$ be the induced
finite covering family of $B'$ in $\cA_{pA'}$, and for
every $i\in I$ let also $k'_i:A'_i\times_{A'}B''\to A'_i$
be the second projection. It suffices to check that
$L_\cA(h'\circ f'_i)\in\cS\times_Ag$ for every $i\in I$,
{\em i.e.} that $L_\cA(k\circ h'\circ f'_i)\in\cS$ for
every such $i$. However, $L_\cA(k\circ h'\circ f'_i)=
L_\cA(h\circ k'\circ f'_i)=L_\cA(h\circ f_i\circ k'_i)$,
whence the contention. To check the local character of
$J^{**}_\cA$, suppose that $\cT$ is another sieve of
$\cA[\Sigma^{-1}_\cA]/A$ such that for every
$(t:A'\to A)\in\cS$ we have $\cT\times_At\in J^{**}_\cA(A')$;
we need to show that $\cT\in J^{**}(A)$. In particular, for
every $i\in I$ there exist a cartesian morphism
$h_i:B_i\to A'_i$ in $\cA$ and a finite covering family
$(g_{i\lambda}:B_{i\lambda}\to B_i~|~\lambda\in\Lambda_i)$ in
$(\cA_{pB_i},J_{pB_i})$ such that
$L_\cA(h_i\circ g_{i\lambda})\in\cT\times_AL_\cA(h\circ f_i)$
for every $\lambda\in\Lambda_i$, {\em i.e.}
$L_\cA(h\circ f_i\circ h_i\circ g_{i\lambda})\in\cT$ for every
such $\lambda$. Now, since $\cC$ is cofiltered, we may find
a morphism $\phi:Y\to pA'$ in $\cC$ such that for every $i\in I$
there exists a morphism $\phi_i:Y\to pB_i$ with
$\phi=p(h_i)\circ\phi_i$. Pick then any cartesian morphism
$k:C\to A'$ with $p(k)=\phi$, and a cartesian morphism
$k_i:C_i\to B_i$ such that $p(k_i)=\phi_i$, for every $i\in I$.
With this notation, for every $i\in I$ there exists a unique
morphism $f'_i:C_i\to C$ in $\cA_Y$ such that
$k\circ f'_i=f_i\circ h_i\circ k_i$. Since $\sc_\phi$ is
continuous for the topologies $J_Y$ and $J_{pA'}$, it follows
easily that $(f'_i:C_i\to C~|~i\in I)$ is a finite covering
family in $(\cA_Y,J_Y)$ (details left to the reader). Lastly,
for every $i\in I$ and $\lambda\in\Lambda_i$ pick a cartesian
morphism $k_{i\lambda}:C_{i\lambda}\to B_{i\lambda}$ with
$p(k_{i\lambda})=\phi_i$.
Then for every such $i$ and $\lambda$ there exists a unique
morphism $g'_{i\lambda}:C_{i\lambda}\to C_i$ in $\cA_Y$ such that
$g_{i\lambda}\circ k_{i\lambda}=k_i\circ g'_{i\lambda}$, and arguing
as in the foregoing, we easily see that
$(g'_{i\lambda}~|~\lambda\in\Lambda_i)$ is a finite covering
family in $(\cA_Y,J_Y)$, for every $i\in I$. We conclude that
$(f'_{i\lambda}:=f'_i\circ g'_{i,\lambda}:C_{i\lambda}\to C~|~
i\in I,\lambda\in\Lambda_i)$ is a finite covering family
in $(\cA_Y,J_Y)$ as well. But we have
$k\circ f'_{i\lambda}=f_i\circ h_i\circ g_{i\lambda}\circ k_{i\lambda}$
for every $i\in I$ and $\lambda\in\Lambda_i$; since $h\circ k$
is a cartesian morphism, the assertion follows. Thus,
$J^{**}_\cA$ is indeed a topology on $\cA[\Sigma_\cA^{-1}]$.
In order to prove that $L_\cA\circ i_X$ is continuous for
$J^{**}_\cA$ and $J_X$, it suffices now to check that for
every covering family
$(f_\lambda:B_\lambda\to B~|~\lambda\in\Lambda)$ in $\cA_X$,
the family $(L_\cA(f_\lambda)~|~\lambda\in\Lambda)$ covers
$B$ for the topology $J^{**}_\cA$ (lemma
\ref{lem_crit-continuity}). Since, as we have already
observed, $B$ is quasi-compact for the topology $J_X$, we
are easily reduced to the case where $\Lambda$ is a finite set
(details left to the reader); but then the assertion follows
immediately from the definition of $J^{**}_\cA(B)$.
\end{pfclaim}

We shall show that the natural morphism of presheaves
$L_{\cA!}F\to(L_{\cA!}F)^+$ is an isomorphism (notation
of \eqref{subsec_asso-topoi}); the lemma will follow
immediately. To this aim, claim \ref{cl_cofinal-sieves}
reduces to checking the following assertion. For every
$A\in\Ob(\cA)$ and every finite covering family
$(f_i:A_i\to A~|~i\in I)$ in $(\cA_{pA},J_{pA})$,
the natural map :
$$
L_{\cA!}FA\to E:=\Equal\Bigl( \prod_{i\in I}
\xymatrix{L_{\cA!}F(A_i)\ar@<-.5ex>[r] \ar@<.5ex>[r] &} \!\!\!\!\!\!\!
\prod_{(i,j)\in I\times I}\!\!\!L_{\cA!}F(A_i\times_AA_j)
\Bigr)
$$
is bijective. Thus, let $\cI$ be the finite category whose set
of objects is the disjoint union of $I$ and $I\times I$, and with
two morphisms $i\leftarrow(i,j)\to j$ for every $(i,j)\in I\times I$
(and of course, the identity morphism of each object); the covering
family $(f_i~|~i\in I)$ induces a functor $A_\bullet:\cI^o\to\cA_{pA}$
given by the rule : $i\mapsto A_i$ and
$(i,j)\mapsto A_{ij}:=A_i\times_AA_j$ for every $i,j\in I$, and
which assigns to every morphism $l:(i,j)\to k$ of $\cI$ the
projection $\pi_l:A_{ij}\to A_k$. We consider the functor
$$
\Psi:\cI\times(\cC/pA)\to\cA
$$
that assigns to every pair $(t,\phi/pA)\in\Ob(\cI)\times\Ob(\cC/pA)$
the object $\sc_\phi(A_t)$, and to every morphism
$(l,\psi):(t,\phi/pA)\to(t',\phi'/pA)$ the composition :
$$
\sc_\phi(A_t)
\xrightarrow{\blambda(\one_{A_t},\psi)}\sc_{\phi'}(A_t)
\xrightarrow{\sc_{\phi'}(\pi_l)}\sc_{\phi'}(A_{t'}).
$$
For every $(t,\phi)\in\Ob(\cI\times(\cC/pA))$, let also $[t]$
(resp. $[\phi]$) be the subcategory of $\cI$ (resp. of
$\cC/pA$) with $\Ob([t]):=\{t\}$ (resp. with
$\Ob([\phi])=\{\phi\}$), and denote by $\Psi_t:\cC/pA\to\cA$
(resp $\Psi^\phi:\cI\to\cA$) the restriction of $\Psi$ to
$[t]\times(\cC/pA)\isom\cC/pA$ (resp. to
$\cI\times[\phi]\isom\cI$). With this notation,
$L_{\cA!}F(A_t)$ represents the colimit of $F\circ\Psi^o_t$,
for every $t\in\Ob(\cI)$, and since $\cI$ is finite and
$(\cC/pA)^o$ is filtered, we get natural identifications :
$$
E=\lim_{t\in\Ob(\cI)}\colim_{(\cC/pA)^o}F\circ\Psi^o_t\isom
\colim_{\phi\in\Ob(\cC/pA)}\lim_\cI F\circ(\Psi^\phi)^o.
$$
But since $\sc_\phi$ is left exact and continuous for the
topologies $J_X$ and $J_Y$, for every morphism $\phi:X\to Y$
of $\cC$, and since the restriction of $F$ to $\cA_X$ is
a sheaf for the topology $J_X$, for every $X\in\Ob(\cC)$,
we see that $F(\sc_\phi A)$ represents the limit of the functor
$F\circ(\Psi^\phi)^o$, for every such $\phi$, so the functor
$(\cC/pA)^o\to\Set$ given by the rule
$\phi\mapsto\lim_\cI F\circ(\Psi^\phi)^o$ is naturally identified
with $F\circ(\st_A\circ\Phi_A)^o$. Summing up, we conclude that
$E$ represents the colimit of $F\circ(\st_A\circ\Phi_A)^o$,
whence (i).

(ii.b): By claim \ref{cl_split-fibrations}, we may assume
that $\cE=\cFib(E_\bullet)$ for a presheaf of categories
$E_\bullet$ on $\cA$. Then recall that $\cF:=\Fib(L_\cA)_!\cE$
is the fibration associated with the strict pseudo-functor
computed by the strong left $2$-Kan extension :
$$
F_\bullet:=2\tdu\!\!\!\int^{L_\cA^o}E_\bullet
\qquad
A\mapsto F_A:=\Pscolim{(A/L_\cA\cA)^o}E_\bullet\circ\st^o_A.
$$
However, since $(A/L_\cA\cA)^o$ is filtered, the $2$-colimit
of $E_\bullet\circ\st^o_A$ is represented by the colimit of
the same functor, for every $A\in\Ob(\cA)$ (example
\ref{ex_filter-2-colim-in-Cat}(iv)), and in view of claim
\ref{cl_Phi-final}, this is also the colimit of the functor
$E_\bullet\circ(\st_A\circ\Phi_A)^o:(\cC/pA)^o\to\bCat$, for every
such $A$.

It suffices to show that the natural functor $\sC(\cF)\to\cF^+$
(see \eqref{subsec_functoriality-of-+}) is an equivalence of
categories, and in view of claim \ref{cl_cofinal-sieves} (and
of corollary \ref{cor_fibrations}(i)) we are reduced to checking
the following assertion. For every $A\in\Ob(\cA)$ and every
finite covering family
$f_\bullet:=(f_i:A_i\to A~|~i\in I)$ in $(\cA_{pA},J_{pA})$,
the natural functor
$$
F_A\isom\cF(A)\to\Desc(\cF,L_\cA(f_\bullet))
$$
is an equivalence. Now, every morphism $\psi/A:\phi\to\phi'$
of $\cC/pA$ induces a functor
$$
\delta_{\psi/A}:
\Desc(\cE,\sc_{\phi'}(f_\bullet))\to\Desc(\cE,\sc_\phi(f_\bullet))
$$
as follows. Notice that, since $\sc_\phi$ is left exact,
a descent datum for $\cE$ relative to the covering
$\sc_\phi(f_\bullet):=(\sc_\phi(f_i)~|~i\in I)$
is a pair $(e_\bullet,\omega_{\bullet\bullet})$, where
$e_\bullet:=(e_i~|~i\in I)$ is a system of objects
$e_i\in\Ob(E_{\sc_\phi A})$ for every $i\in I$, and
$\omega_{\bullet\bullet}$ is a system of isomorphisms
$\omega_{ij}:E_{\sc_\phi(\pi_{ij})}(e_i)\isom E_{\sc_\phi(\pi_{ji})}(e_j)$
fulfilling the usual cocycle condition (here we denote by
$\pi_{ij}:A_{ij}\to A_i$ and $\pi_{ji}:A_{ij}\to A_j$ the
projections). A morphism $h_\bullet:
(e_\bullet,\omega_{\bullet\bullet})\to(e_\bullet,\omega_{\bullet\bullet})$
of such descent data is a system of morphisms
$(h_i:e_i\to e'_i~|~i\in I)$ compatible in the obvious
fashion with $\omega_{ij}$ and $\omega'_{ij}$. Then
$\delta_{\psi/A}$ is given by the rules :
$$
\begin{aligned}
(e_\bullet,\omega_{\bullet\bullet})&\,\mapsto
((E_{\blambda(\one_{A_i},\psi)}(e_i)~|~i\in I),
(E_{\blambda(\one_{A_{ij}},\psi)}(\omega_{ij})~|~i,j\in I)) \\
h_\bullet&\,\mapsto((E_{\blambda(\one_{A_i},\psi)}(h_i)~|~i\in I)
\end{aligned}
$$
for every object $(e_\bullet,\omega_{\bullet\bullet})$ and every
morphism $h_\bullet$ of $\Desc(\cE,\sc_{\phi'}(f_\bullet))$.
We thus obtain a functor $\delta:(\cC/pA)^o\to\bCat$
that assigns to every $\phi\in\Ob(\cC/pA)$ the category
$\Desc(\cE,\sc_\phi(f_\bullet))$, and to every morphism $\psi/A$
the functor $\delta_{\psi/A}$, and a direct inspection of the
definitions yields a natural identification
$$
\Desc(\cF,L_\cA(f_\bullet))\isom\colim_{(\cC/pA)^o}\delta.
$$
Lastly, since $\cE$ is a stack on $(\cA[\Sigma^{-1}_\cA],J^*_\cA)$,
the natural functor $E_{\sc_\phi A}\to\Desc(\cE,\sc_\phi(f_\bullet))$
is an equivalence for every $\phi\in\Ob(\cC/pA)$, and we get
an induced equivalence
$$
\Pscolim{(\cC/pA)^o}E_\bullet\circ(\st\circ\Phi_A)^o\isom
\Pscolim{(\cC/pA)^o}\delta.
$$
But again, these $2$-colimits are represented by the colimits
of the same functors, whence the contention.
\end{proof}

\sset\subsubsection{}\label{subsec_bc-for-limit-site}
Keep the notation of \eqref{subsec_localization-morph-of-sites},
and let now $\phi:(\cA,p,J_\bullet)\to(\cA',p',J'_\bullet)$ be a
morphism of small fibred lex-sites over the small category $\cC$.
We have $\phi(\Sigma_{\cA'})\subset\Sigma_\cA$, whence a commutative
diagram of categories, for every $X\in\Ob(\cC)$ :
$$
\xymatrix{ \cA'_X \ar[r]^-{i'_X} \ar[d]_{\phi_X} &
\cA' \ar[r]^-{L_{\cA'}} \ar[d]_\phi &
\cA'[\Sigma^{-1}_{\cA'}] \ar[d]^{\phi_\Sigma} \\
\cA_X \ar[r]^-{i_X} & \cA \ar[r]^-{L_\cA} & \cA[\Sigma^{-1}_\cA]
}$$
where $i_X$ and $i'_X$ are the inclusion functors. Let as well
$\blambda'$ be a unital cleavage for $p'$, and $\sc'$ its
associated pseudo-functor, which factors likewise through
a pseudo-functor $\tilde\sc{}':\cC\to\lex.\Site$; then $\phi$
corresponds to a unique pseudo-natural transformation
$\omega:\sc'\Rightarrow\sc$, and the latter in turn can also
be regarded as a $2$-cell
$\tilde\omega:\tilde\sc\Rightarrow\tilde\sc{}'$ of $\lex.\Site$.
Now, since $\cA[\Sigma^{-1}_\cA]$ is a strong $2$-colimit for
$\sc$, there exists a unique functor
$\psi:\cA'[\Sigma^{-1}_{\cA'}]\to\cA[\Sigma^{-1}_\cA]$ such that
$\tau\odot\sF_\psi=\omega$; clearly we must then have
$\psi=\phi_\Sigma$. On the other hand, since
$\cA[\Sigma^{-1}_\cA]$ is a strong $2$-limit for $\tilde\sc$,
there exists a unique morphism of sites
$\tilde\psi:(\cA[\Sigma^{-1}_\cA],J^*_\cA)\to
(\cA'[\Sigma^{-1}_{\cA'}],J^*_{\cA'})$ such that
$\sF_{\tilde\psi}\odot\tilde\tau=\tilde\omega$; thus, we must
have $\tilde\psi=\psi$, and this shows that
$\phi_\Sigma$ is a morphism of sites as well,
so we get a commutative diagram of morphisms of sites :
$$
\xymatrix@C+20pt{ (\cA[\Sigma^{-1}_\cA],J^*_\cA)
\ar[r]^-{L_\cA} \ar[d]_{\phi_\Sigma} &
(\cA,J_\cA) \ar[d]^\phi \\
(\cA'[\Sigma^{-1}_{\cA'}],J^*_{\cA'}) \ar[r]^-{L_{\cA'}} &
(\cA',J_{\cA'})
}$$
which we orient by the identity transformation $\one_{\phi\circ L_\cA}$.
Following \eqref{subsec_rotate-the-cube}, the resulting oriented
diagram then induces an oriented square of $2$-categories :
$$
\xymatrix@C+30pt{ \Stack(\cA[\Sigma^{-1}_\cA],J^*_\cA)
\drtwocell\omit{_\ \ \ \ \ \ \ \ \ \ \ \ \ \ \St(\one_{\phi\circ L_\cA})^\gamma_*}
\ar[r]^-{\St(L_\cA)_*} \ar[d]_{\St(\phi_\Sigma)_*} &
\Stack(\cA,J_\cA) \ar[d]^{\St(\phi)_*} \\
\Stack(\cA'[\Sigma^{-1}_{\cA'}],J^*_{\cA'}) \ar[r]_-{\St(L_{\cA'})_*} &
\Stack(\cA',J_{\cA'})
}$$
which in turn can be regarded as an oriented square of links
of the $2$-category $\sV\tdu\overline{2\tdu\bCat}$, for a
suitable universe $\sV$ : see \eqref{subsec_return-to-wlinks}.

\begin{proposition}\label{prop_bc-lim-and-tot}
In the situation of \eqref{subsec_bc-for-limit-site}, suppose
that every object of $\cA$ is quasi-compact for the topology
$J_\cA$. Then $\Upsilon(\St(\one_{\phi\circ L_\cA})^\gamma_*)$ is a
pseudo-natural equivalence.
\end{proposition}
\begin{proof} Let $\cE$ be any stack on $(\cA,J_\cA)$; as usual,
we may assume that $\cE=\cFib(E_\bullet)$ for a presheaf of
categories $E_\bullet$ on $\cA$, and by virtue of lemma
\ref{lem_stacks-and-sheaves-of-cats}(ii) we may also assume
that $E_\bullet$ is a sheaf of categories on $(\cA,J_\cA)$.
According to corollary \ref{cor_conditions-for-trivial-bc}(iv),
and the discussion of \eqref{subsec_from-beta-dagger-to-tilde},
we are then reduced to checking :

\begin{claim}(i)\ \ $\cFib((L^\sim_\cA,\bCat)^*E_\bullet)$
is a stack on $(\cA[\Sigma^{-1}_\cA],J^*_\cA)$.

(ii)\ \
$\Upsilon((\one_{\phi\circ L_\cA})^\sim_*)$ is an isomorphism
of functors.
\end{claim}
\begin{pfclaim}[](i): On one hand the natural morphism
$(L_\cA,\bCat)_!E_\bullet\to(L^\sim_\cA,\bCat)^*E_\bullet$ is an
isomorphism, by proposition \ref{prop_local-morph-of-sites}(ii.a);
on the other hand, from claim \ref{cl_Phi-final} and lemma
\ref{lem_filtered-final}(ii) we see that the category
$A/L_\cA\cA$ is cofiltered for every $A\in\Ob(\cA)$, hence
remark \ref{rem_alternative-inverse-image}(v) yields an
equivalence of categories
$\cFib((L_\cA,\bCat)_!E_\bullet)\isom\Fib(L_\cA)_!(\cE)$.
Lastly, $\Fib(L_\cA)_!(\cE)$ is a stack, by proposition
\ref{prop_local-morph-of-sites}(ii.b), whence the assertion.

(ii): By definition, we have :
$$
\Upsilon((\one_{\phi\circ L_\cA})^\sim_*)=
(\eps^{L_{\cA'}}*\tilde\phi_{\Sigma*}\tilde L{}^*_\cA)\odot
(\tilde L{}_{\cA'}^*\tilde\phi_**\eta^{L_\cA})
$$
where $\eta^{L_\cA}$ (resp. $\eps^{L_{\cA'}}$) is the unit
(resp. the counit) of an adjunction for the pair
$(\tilde L{}^*_\cA,\tilde L_{\cA*})$ (resp.
$(\tilde L{}^*_{\cA'},\tilde L_{\cA'*})$). However,
$\tilde L_{\cA'*}$ is fully faithful (proposition
\ref{prop_local-morph-of-sites}(i)), hence $\eps^{L_{\cA'}}$
is an isomorphism (proposition \ref{prop_fullfaith-adjts}(iii)).
Thus, we are reduced to checking that
$\tilde L{}_{\cA'}^*\tilde\phi_**\eta^{L_\cA}$ is an isomorphism.
However, since also $\eps^{L_\cA}$ is an isomorphism, we know
already that $\tilde L{}^*_\cA*\eta^{L_\cA}$ is an isomorphism,
due to the triangular identities of \eqref{subsec_adj-pair}.
Thus, we are further reduced to exhibiting an {\em arbitrary}
isomorphism of functors :
$$
\tilde L{}_{\cA'}^*\circ\tilde\phi_*\isom
\tilde\phi_{\Sigma*}\circ\tilde L{}^*_\cA.
$$
Hence, let $F$ be any sheaf on $(\cA,J_\cA)$, and $A'\in\Ob(\cA')$;
from the proof of proposition \ref{prop_local-morph-of-sites}(i)
we see that $\tilde L{}_{\cA'}^*\circ\tilde\phi_*F(A')$ represents
the colimit of the functor
$$
\Psi:(\cC/p'A')^o\to\Set
\qquad
(X\xrightarrow{f}p'A')\mapsto F\phi(\sc'_fA')
\qquad
(f\xrightarrow{g/p'A'}f')\mapsto F\phi(\blambda'(\one_{A'},g))
$$
whereas $\tilde\phi_{\Sigma*}\circ\tilde L{}^*_\cA F(A')$
represents the colimit of the functor
$$
\Psi':(\cC/p'A')^o\to\Set
\qquad
(X\xrightarrow{f}p'A')\mapsto F\sc_f(\phi A')
\qquad
(f\xrightarrow{g/p'A'}f')\mapsto F(\blambda'(\one_{\phi A'},g)).
$$
Now, let $\tau^\omega$ denote the coherence constraint of
the pseudo-natural transformation $\omega$ associated with
$\phi$ as in \eqref{subsec_bc-for-limit-site}, so that
$\omega_X:\cA_X\to\cA'_X$ is the restriction of $\phi$,
for every $X\in\Ob(\cC)$; we claim that the following system
of maps yields an isomorphism of functors $\Psi\isom\Psi'$ :
\set\begin{equation}\label{eq_tired-knee}
(F(\tau^\omega_{f,A'}):F\phi(\sc'_f A')\isom F\sc_f(\phi A')~|~
f\in\Ob(\cC/p'A')).
\end{equation}
More precisely, we show that for every two objects $f:X\to p'A'$
and $f':X'\to p'A'$ and every morphism $g/p'A':f\to f'$
of $\cC/p'A'$ the following diagram commutes :
\set\begin{equation}\label{eq_milk-this-later}
{\spreaddiagramcolumns{+20pt}\diagram
\sc_f(\phi A') \ar[r]^-{\tau^\omega_{f,A'}}
\ar[d]_{\blambda(\one_{\phi A'},g/p'A')} &
\phi(\sc'_f A') \ar[d]^{\phi(\blambda'(\one_{A'},g/p'A'))} \\
\sc_{f'}(\phi A') \ar[r]^-{\tau^\omega_{f',A'}} & \phi(\sc'_{f'}A').
\enddiagram}\end{equation}
To this aim, denote by $\gamma^{\sc'}$ the coherence constraint
of $\sc'$, and recall that, under the natural identification
$\cA'\isom\cFib(\sc')$, the morphism $\blambda'(\one_{A'},g/p'A')$
of $\cA'$ corresponds to the morphism
$$
(g,\gamma^{\sc'\,-1}_{(g,f'),A'}:\sc'_fA'\to\sc'_g\sc'_{f'}A'):
(X,\sc'_fA')\to(X',\sc'_{f'}A').
$$
Likewise we describe $\blambda(\one_{\phi A'},g)$ in terms of the
coherence constraint $\gamma^\sc$ of $\sc$; since $\blambda$ and
$\blambda'$ are unital, the assertion comes down to the
commutativity of the diagram :
$$
\xymatrix@C+40pt{ \sc_g\sc_{f'}(\phi A')
\ar[rr]^-{\gamma^\sc_{(g,f'),\phi A'}}
\ar[d]_{\sc_g(\tau^\omega_{f',A'})} & &
\sc_f(\phi A') \ar[d]^{\tau^\omega_{f,A'}} \\
\sc_g\circ\phi(\sc'_{f'}A') \ar[r]^-{\tau^\omega_{g,\sc'_{f'}A'}} &
\ar[r]^-{\omega_X(\gamma^{\sc'}_{(g,f'),A'})} \phi(\sc'_g\sc'_{f'}A')
& \phi(\sc'_fA')
}$$
which in turn follows from the coherence axiom for $\tau^\omega$.
So finally, the colimit of the system \eqref{eq_tired-knee}
yields a bijection $t_{F,A'}:\tilde L{}_{\cA'}^*\circ\tilde\phi_*F(A')
\isom\tilde\phi_{\Sigma*}\circ\tilde L{}^*_\cA F(A')$ for every
$A'\in\Ob(\cA')$, and we need to check that the rule :
$A'\mapsto t_{F,A'}$ defines a morphism of sheaves $t_F$. Thus, let
$$
(\beta^A_{F,f}:F(\sc_fA)\to\tilde L{}^*_\cA F(A)~|~f\in\Ob(\cC/pA))
\qquad
\text{for every $A\in\Ob(\cA)$}
$$
be a fixed system of universal cocones, and $h:B\to A$ any
morphism in $\cA[\Sigma^{-1}_\cA]$; for a given object
$f:X\to pA$ of $\cC/pA$, we may find a morphism
$g:Y\to pB$ in $\cC$ and a morphism $h_{g,f}:\sc_gB\to\sc_fA$
in $\cA$ such that $L_\cA(\blambda(A,f))^{-1}\circ h=
L_\cA(h_{g,f})\circ L_\cA(\blambda(B,g))^{-1}$, and by unraveling
the definitions, it is easily seen that :
$$
\tilde L{}^*_\cA F(h)\circ\beta^A_{F,f}=\beta^B_{F,g}\circ Fh_{g,f}.
$$
Likewise, fix a system of universal cocones :
$$
(\beta'^{A'}_{F,f}:F\phi(\sc'_fA')\to
\tilde L{}^*_{\cA'}\tilde\phi_*F(A')~|~f\in\Ob(\cC/p'A'))
\qquad
\text{for every $A'\in\Ob(\cA')$}.
$$
If $h':B'\to A'$ is any morphism of $\cA'[\Sigma^{-1}_{\cA'}]$,
and $f:X\to p'A'$ any object of $\cC/p'A'$, we may find
a morphism $g:Y\to p'B'$ in $\cC$ and a morphism
$h'_{g,f}:\sc'_gB'\to\sc'_fA'$ of $\cA'$ such that
$L_{\cA'}(\blambda'(A,f))^{-1}\circ h'=
L_{\cA'}(h'_{g,f})\circ L_{\cA'}(\blambda'(B,g))^{-1}$, and again
it is easily seen that :
$$
\tilde L{}^*_{\cA'}\tilde\phi_*F(h')\circ\beta'^{A'}_{F,f}=
\beta'^{B'}_{F,g}\circ F\phi(h'_{g,f}).
$$
Now, for every such $h'$ we need to check the identity :
\set\begin{equation}\label{eq_natucha}
t_{F,B'}\circ\tilde L{}_{\cA'}^*\tilde\phi_*F(h')=
\tilde\phi_{\Sigma*}\tilde L{}^*_\cA F(h')\circ t_{F,A'}
\end{equation}
and it suffices to check that the sides of \eqref{eq_natucha}
agree after composition with $\beta'^{A'}_{F,f}$, for every
$f\in\Ob(\cC/p'A')$. Notice that if $g:Y\to p'B'$ is a
suitable choice for constructing $h'_{g,f}$, then the same
choice can be used for constructing $(\phi_\Sigma h')_{g,f}$
(details left to the reader). We compute:
$$
\begin{aligned}
t_{F,B'}\circ\tilde L{}_{\cA'}^*\tilde\phi_*F(h')\circ\beta'^{A'}_{F,f}
&\,=t_{F,B'}\circ\beta'^{B'}_{F,g}\circ F\phi(h'_{g,f})=
\beta^{\phi B'}_{F,g}\circ F(\tau^\omega_{g,B'})\circ F\phi(h'_{g,f}) \\
\tilde\phi_{\Sigma*}\tilde L{}^*_\cA F(h')\circ t_{F,A'}\circ\beta'^{A'}_{F,f}
&\,=\tilde\phi_{\Sigma*}\tilde L{}^*_\cA F(h')\circ
\beta^{\phi A'}_{F,f}\circ F(\tau^\omega_{f,A'})=
\beta^{\phi B'}_{F,g}\circ F(\phi_\Sigma h')_{g,f}\circ F(\tau^\omega_{f,A'}).
\end{aligned}
$$
Thus, we are reduced to checking that we may choose
$(\phi_\Sigma h')_{g,f}$ such that :
$$
\phi(h'_{g,f})\circ\tau^\omega_{g,B'}=
\tau^\omega_{f,A'}\circ(\phi_\Sigma h')_{g,f}.
$$
In other words, we have to check the identity :
\set\begin{equation}\label{eq_pio-pio}
L_\cA((\tau^\omega_{f,A'})^{-1}\circ\phi(h'_{g,f})\circ\tau^\omega_{g,B'})
\circ L_\cA(\blambda(\phi B',g))^{-1}=
L_\cA(\blambda(\phi A',f))^{-1}\circ(\phi_\Sigma h').
\end{equation}
However, if we regard $g$ as a morphism $g/p'B':g\to\one_{p'B'}$
of $\cC/p'B'$, we get
$$
\blambda'(\one_{B'},g/p'B')=\blambda'(B',g)
\qquad\text{and}\qquad
\blambda(\one_{\phi B'},g/p'B')=\blambda(\phi B',g).
$$
Likewise, we have corresponding identities for $\blambda'(A',g)$
and $\blambda(\phi A',g)$, by regarding $f$ as a morphism
$f/p'A':f\to\one_{p'A'}$ of $\cC/p'A'$. Then
\eqref{eq_milk-this-later} yields the identities
$$
\phi(\blambda'(A',f))\circ\tau^\omega_{f,A'}=\blambda(\phi A',f)
\qquad
\phi(\blambda'(B',g))\circ\tau^\omega_{g,B'}=\blambda(\phi B',g).
$$
(notice that $\tau^\omega_{\one_{A'},A'}=\one_{\phi A'}$ and
$\tau^\omega_{\one_{B'},B'}=\one_{\phi B'}$, since $\sc$ and
$\sc'$ are unital : see remark \ref{rem_unital}(ii)).
Summing up, we deduce that \eqref{eq_pio-pio} is equivalent
to the identity :
$$
L_\cA(\phi(\blambda'(A',f))\circ\phi(h'_{g,f}))
\circ L_\cA(\phi(\blambda'(B',g)))^{-1}=\phi_\Sigma h'
$$
which follows immediately from the definition of $h'_{g,f}$
and the identity $\phi_\Sigma\circ L_{\cA'}=L_\cA\circ\phi$.

This conclude the construction of the morphism of sheaves
$t_F:\tilde L{}_{\cA'}^*\circ\tilde\phi_*F
\isom\tilde\phi_{\Sigma*}\circ\tilde L{}^*_\cA F$. To conclude
the proof, it remains to check that the rule $F\mapsto t_F$
yields the sought isomorphism of functors. Thus, let
$\nu:F\to G$ be a morphism of sheaves on $(\cA,J_\cA)$;
we come down to checking the identity :
\set\begin{equation}\label{eq_shark}
\tilde\phi_{\Sigma*}\tilde L{}^*_\cA(\nu)_{A'}\circ t_{F,A'}=
t_{G,A'}\circ\tilde L{}^*_\cA\tilde\phi_*(\nu)
\qquad
\text{for every $A'\in\Ob(\cA')$}
\end{equation}
and again, it suffices to check that the two sides of
\eqref{eq_shark} agree after composition with $\beta'^{A'}_{F,f}$,
for every $f\in\Ob(\cC/p'A')$. We compute :
$$
\begin{aligned}
\tilde\phi_{\Sigma*}\tilde L{}^*_\cA(\nu)_{A'}\circ t_{F,A'}
\circ\beta'^{A'}_{F,f}&\,=\tilde\phi_{\Sigma*}\tilde L{}^*_\cA(\nu)_{A'}
\circ\beta^{\phi A'}_{F,f}\circ F(\tau^\omega_{f,A'})=
\beta^{\phi A'}_{G,f}\circ\nu_{\sc_f(\phi A')}\circ F(\tau^\omega_{f,A'}) \\
t_{G,A'}\circ\tilde L{}^*_\cA\tilde\phi_*(\nu)\circ\beta'^{A'}_{F,f}
&\,=t_{G,A'}\circ\beta'^{A'}_{G,f}\circ\nu_{\phi(\sc'_fA')}=
\beta^{\phi A'}_{G,f}\circ G(\tau^\omega_{f,A'})\circ\nu_{\phi(\sc'_fA')}.
\end{aligned}
$$
So we come down to showing that
$\nu_{\sc_f(\phi A')}\circ F(\tau^\omega_{f,A'})=
G(\tau^\omega_{f,A'})\circ\nu_{\phi(\sc'_fA')}$. The latter is clear,
since $\nu$ is a morphism of sheaves.
\end{pfclaim}
\end{proof}

\sset\subsubsection{Stacks on fibred topoi}
\label{subsec_one-more-reduction}
Recall that the pseudo-functor $\sT$ of \eqref{subsec_T-and-Can}
assigns to every site $C$ a topos $\sT(C)$ with an isomorphism
$\omega_C:C^\sim\isom\sT(C)$, and the unit of the $2$-adjoint
pair $(\Can,\sT)$ exhibited in the proof of theorem
\ref{th_adj-topos-site} assigns to every topos $T$ an
equivalence $\shh_{T*}:\isom\sT\circ\Can(T)$ whose composition
with $\omega_{\Can(T)}^{-1}$ equals the Yoneda imbedding
$h_T:T\isom T^\sim$. Also, for every natural transformation
$\beta:u\Rightarrow v$ of morphisms of topoi $u,v:T'\to T$,
the isomorphisms $\omega_{\Can(T)}$ and $\omega_{\Can(T')}$
identify $\sT\circ\Can(\beta)$ with $(\beta^\dagger)^\sim_*$,
where $\beta^\dagger$ is the adjoint transformation to $\beta$.
Now, consider any oriented square of morphisms of topoi :
$$
\xymatrix{ T' \ar[r]^-u \ar[d]_{f'}
\drtwocell\omit{_\ \beta} & T \ar[d]^f \\
S' \ar[r]_-v & S
}$$
({\em i.e.} $\beta$ is a natural transformation
$f_*\circ u_*\Rightarrow v_*\circ f'_*$). There follows
a diagram of oriented squares of topoi :
$$
\xymatrix{ (\Can\ T')^\sim \ar[rrr]^-{(u^*)^\sim}
\ar[ddd]_{(f'^*)^\sim} & &
\dltwocell\omit{^\ } & (\Can\ T)^\sim \ar[ddd]^{(f^*)^\sim} \\
& \dltwocell\omit{^\ } T' \ar[r]^-u \ar[d]_{f'} \ar[lu]^{h_{T'}}
\drtwocell\omit{_\ \beta} & T \ar[d]^f \ar[ru]_{h_T}
\drtwocell\omit{_\ } \\
& S' \ar[r]_-v \ar[ld]_{h_{S'}}\drtwocell\omit{^\ }
& S \ar[rd]^{h_S} & \\
(\Can\ S')^\sim \ar[rrr]_-{(v^*)^\sim} & & & (\Can\ S)^\sim
}$$
whose unmarked orientations are identified, via the
isomorphisms $\omega_\bullet$, with the coherence constraints
of the pseudo-natural transformation $\shh_{\bullet*}$. We
complete it by inserting the orientation $(\beta^\dagger)^\sim_*$
for the external square. Then the resulting cubical diagram
commutes on $2$-cells, in the sense of
\eqref{subsec_transfer-base-change}. Since the four diagonal
arrows are equivalences (theorem \ref{th_canon-topos}(iv)) it
follows that $\Upsilon(\beta)$ is an isomorphism of functors if
and only if the same holds for $\Upsilon((\beta^\dagger)^\sim_*)$
(remark \ref{rem_transit-base-change}(ii)).

\sset\subsubsection{}
\label{subsec_trivial-base-change-stacks}
Let $\omega:T\to S$ be a morphism of fibred topoi over a
small category $I$; let also $t\in\Ob(I)$ be any object,
and define the functor $\rho^t:\bone\to I$ as in remark
\ref{rem_fibrewise-lims-in-tot-Top}(ii). By direct inspection,
we see that :
$$
\Psi^*_t:=\sTop(T/\rho^t)^*\circ\sTop(\omega)^*=
\omega_t^*\circ\sTop(S/\rho^t)^*.
$$
Indeed, both functors attach to every object
$E_\bullet:I\to\cFib(S^*)$ of $\sTop(S)$ the object
$\omega^*_t(E_t)\in\Ob(T_t)$, and to every morphism
$\beta_\bullet:E_\bullet\to E'_\bullet$, the morphism
$\omega^*_t(\beta_t)$. Similarly, we get :
$$
\Phi_t:=\sTop(S/\rho^t)^*\circ\sTop(\omega)_*=
\omega_{t*}\circ\sTop(T/\rho^t)^*.
$$
The adjoint of $\one_{\Psi^*_t}$ yields therefore an
orientation for the following diagram of topoi :
$$
\xymatrix{ T_t \ar[rr]^-{\sTop(T/\rho^t)} \ar[d]_{\omega_t}
\drrtwocell\omit{_\ \ \ \ \ \one^\dagger_{\Psi^*_t}} & &
\sTop(T) \ar[d]^{\sTop(\omega)} \\
S_t \ar[rr]_-{\sTop(S/\rho^t)} & & \sTop(S)
}$$
which we regard as an oriented square of links in the $2$-category
of categories (see \eqref{subsec_base-change-map}).

\begin{lemma}\label{lem_trivial-one}
With the notation of \eqref{subsec_trivial-base-change-stacks},
we have : $\Upsilon(\one^\dagger_{\Psi^*_t})=\one_{\Phi_t}$.
\end{lemma}
\begin{proof} We may as usual assume that $T$ and $S$ are
unital; then we may replace $\sTop(S)$ by the isomorphic
category $\sTop(S)_*$, and $\sTop(\omega)$ by the adjunction
$\Omega:=(\Omega^*,\Omega_*,\eta^\Omega)$ described in the proof
of proposition \ref{prop_top-adjoint-pair}. Then
$\sTop(S/\rho^t)^*$ is replaced by $\Sigma(\cFib(S_*)/\rho^{to})^*$
(see the proof of corollary \ref{cor_left-right-adj-Top-rho}),
and we come down to considering the oriented square of links :
$$
{\spreaddiagramcolumns{+50pt}\diagram
T_t \ar[r]^-{\sTop(T/\rho^t)} \ar[d]_{\omega_t}
\drtwocell\omit{_\ \ \ \ \ \one^\dagger_{\Lambda^*_t}} &
\sTop(T) \ar[d]^{\Omega} \\
S_t \ar[r]_-{\Sigma(\cFib(S_*)/\rho^{to})} & \sTop(S)_*
\enddiagram}
\qquad
\text{where $\Lambda^*_t:=\Top(T/\rho^t)^*\circ\Omega^*$}
$$
whose bottom horizontal arrow is given by
$\Sigma(\cFib(S_*)/\rho^{to})^*$, its right adjoint, and a
unit for this adjoint pair. Set
$\Phi'_t:=\Sigma(\cFib(S_*)/\rho^{to})^*\circ\Omega_*$;
by proposition \ref{prop_opp-links-and-base-ch}, it suffices
to show that $\Upsilon({}^o\one_{\Lambda^*_t})={}^o\one_{\Phi'_t}$.
However, let $(\eta^{\omega_i},\eps^{\omega_i})$ be the unit and
counit of the adjunction for the pair $(\omega_i^*,\omega_{i*})$
defining the morphism of topoi $\omega_i$, for every
$i\in\Ob(I)$; by definition we have :
$$
\Upsilon({}^o\one_{\Lambda^*_t})=
(\eps^{{}^o\omega_t}*{}^o\Sigma(\cFib(S_*)/\rho^{to})^**{}^o\Omega_*)
\odot({}^o\omega_{t*}*{}^o\sTop(T/\rho^t)^**\eta^{{}^o\Omega})
$$
and recall that $\eps^{{}^o\omega_t}={}^o\eta^{\omega_t}$; likewise,
$\eta^{{}^o\Omega}={}^o\eps^\Omega$, where $\eps^\Omega$ denotes the
counit for the adjoint pair $\Omega$. Thus, we get :
$$
{}^o\Upsilon({}^o\one_{\Lambda^*_t})=
(\omega_{t*}*\sTop(T/\rho^t)^**\eps^\Omega)\odot
(\eta^{\omega_t}*\Sigma(\cFib(S_*)/\rho^{to})^**\Omega_*)=
(\Phi_t*\eps^\Omega)\odot(\eta^{\omega_t}*\Phi_t).
$$
But by inspecting the proof of proposition \ref{prop_top-adjoint-pair}
we see that $\eps^\Omega$ is the natural transformation that assigns
to every object $E_\bullet:I\to\cFib(T^*)$ of $\sTop(T)$ the system
of morphisms
$(\eps^{\omega_i}:\omega^*_i\circ\omega_{i*}E_i\to E_i~|~i\in\Ob(I))$.
Thus, finally, the natural transformation
$(\Phi_t*\eps^\Omega)\odot(\eta^{\omega_t}*\Phi_t)$ attaches to every
such $E_\bullet$ the morphism
$(\omega_{t*}*\eps^{\omega_t})\odot(\eta^{\omega_t}*\omega_{t*})$, which
is indeed $\one_{\omega_{t*}E_t}$, by the triangular identities of
\eqref{subsec_adj-pair}.
\end{proof}

\sset\subsubsection{}\label{subsec_flowers}
Keep the situation of \eqref{subsec_trivial-base-change-stacks},
and to ease notation set
$$
v_t:=\sTop(S/\rho^t)^*
\qquad
v'_t:=\sTop(T/\rho^t)^*
\qquad
\sS:=\Can\circ\sTop:\sPsFun(I,\Topos)\to\Site.
$$
Arguing as in \eqref{subsec_rotate-the-cube}, we deduce an
essentially commutative diagram :
$$
\xymatrix@C+30pt{ \Stack(\Can\,T_t) \ar[d]_{\St(\omega^*_t)_*}
\ar[r]^-{\St(v'_t)_*}
\drtwocell\omit{_\ \ \ \ \ \ \ \ \ \ \ \ \St(\one_{\Psi^*_t})^\gamma_*}
& \Stack(\sS(T)) \ar[d]^{\St(\sTop(\omega)^*)_*} \\
\Stack(\Can\,S_t) \ar[r]_-{\St(v_t)_*} & \Stack(\sS(S))
}$$
which we regard as an oriented square in
$\sLink(\sV\tdu\overline{2\tdu\bCat})$ (see
\eqref{subsec_return-to-wlinks}); then we may state :

\begin{proposition}\label{prop_nerd}
$\Upsilon(\St(\one_{\Psi^*_t})^\gamma_*)$ is a pseudo-natural equivalence.
\end{proposition}
\begin{proof} Let $\cE$ be any stack on the canonical site of
$\sTop(T)$; we need to check that
$\Upsilon(\St(\one_{\Psi^*_t})^\gamma_*)_\cE$ is an equivalence of
categories. To this aim, in view of claim
\ref{cl_split-fibrations} we may replace $\cE$ by the split
fibration $\sC(\cE)$, so that $\cE=\cFib(\cA_\bullet)$ for a
presheaf of categories $\cA_\bullet$ on $\sTop(T)$; then by lemma
\ref{lem_stacks-and-sheaves-of-cats}(ii) we may as well replace
$\cA_\bullet$ by $\cA_\bullet^a$, and assume that $\cA_\bullet$ is
a sheaf of categories on $\sTop(T)$. 

\begin{claim}\label{cl_trivial-one}
$\Upsilon((\one_{\Psi^*_t}^\sim,\bCat)_*)$ is a pseudo-natural
equivalence.
\end{claim}
\begin{pfclaim} As explained in
\eqref{subsec_from-beta-dagger-to-tilde}, it suffices to check
that $\Upsilon((\one_{\Psi^*_t})^\sim_*)$ is an isomorphism of
functors. By \eqref{subsec_one-more-reduction}, the latter holds
if and only if $\Upsilon(\one_{\Psi^*_t}^\dagger)$ is an isomorphism
of functors. This in turn is clear from lemma \ref{lem_trivial-one}.
\end{pfclaim}

In light of claim \ref{cl_trivial-one} and corollary
\ref{cor_conditions-for-trivial-bc}(iv), it suffices
to check that $\cFib((\tilde v{}'_t,\bCat)^*\cA_\bullet)$ is a
stack on the canonical site of $T_t$. For every site $C$, let
$$
\omega_C:\bCat^*(C^\sim)\isom(C,\bCat)^\sim
$$
be the strict and strong $2$-equivalence of
\eqref{eq_upgrade-Cats-and-stars}. We may assume that
$\cA_\bullet=\omega_{\sTop(T)}(\cA^*)$ for a category object
$\cA^*$ of $\Can(\sTop(T))^\sim$ and then we need to check that
$\cFib(\omega_{T_t}\circ\bCat^*(\tilde v{}'^*_t)(\cA^*))$ is a
stack on $\Can(T_t)$. Since $h_{\sTop(T)}:\sTop(T)\to\Can(\sTop(T))^\sim$
is an equivalence, we may assume that
$\cA^*=\bCat^*(h_{\sTop(T)})(\cB^*)$ for some
$\cB^*\in\Ob(\bCat^*(\sTop(T)))$; by corollary
\ref{cor_two-U-sites}(i.b), it then suffices to check that
$\cFib(\omega_{T_t}\circ\bCat^*(h_{T_t})\circ\bCat^*(v'_t)(\cB^*))$
is a stack. However, theorem \ref{th_top-is-a-topos}(i) yields
an equivalence of categories $a:\sTop(T)\isom(\Can(T),J)^\sim$
and a morphism of sites $b:\Can(\sTop(T))\to(\Can(T),J)$
with isomorphisms of functors :
$$
h_{T_t}\circ v'_t\isom\tilde\imath_{t*}\circ a
\qquad
a\isom\tilde b_*\circ h_{\sTop(T)}
$$
where $i_t:T_t\to\Can(T)$ is the inclusion functor.
There follow isomorphisms of fibrations :
$$
\begin{aligned}
\cFib(\omega_{T_t}\circ\bCat^*(h_{T_t})\circ\bCat^*(v'_t)\cB^*)
\isom&\,\cFib(\omega_{T_t}\circ\bCat^*(\tilde\imath_{t*})\circ
\bCat^*(a)\cB^*) \\
\isom&\,
\cFib((\tilde\imath_t,\bCat)_*\circ\omega_{(\Can(T),J)}\circ
\bCat^*(a)\cB^*) \\
\isom&\,
\Fib(i_t)^*(\cFib(\omega_{(\Can(T),J)}\circ\bCat^*(a)\cB^*))
\end{aligned}
$$
and by corollary \ref{cor_lims-and-totsites}, it then suffices
to check that $\cFib(\omega_{(\Can(T),J)}\circ\bCat^*(a)\cB^*)$
is a stack on $(\Can(T),J)$. We are then further reduced
to showing that $\cFib((\tilde b,\bCat)_*\cA_\bullet)$ is
a stack on $(\Can(T),J)$, or equivalently, that the same
holds for $\Fib(b)^*(\cFib(\cA_\bullet))=\Fib(b)^*(\cE)$.
The latter holds, by corollary \ref{cor_another-condition}.
\end{proof}

\sset\subsubsection{}\label{subsec_first-for-topoi}
Consider now an oriented diagram of fibred topoi over a small
category $I$ :
$$
\xymatrix{ T' \ar[r]^-{\mu} \ar[d]_{\omega'}
\drtwocell\omit{_\ \ \Xi} & T \ar[d]^\omega \\
S' \ar[r]_\nu & S
}$$
{\em i.e.} $\Xi:\omega\circ\mu\leadsto\nu\circ\omega'$
is a modification. For every $t\in\Ob(I)$ we get an oriented
diagram :
$$
\xymatrix@C+50pt{ S_t \ar[rrr]^-{\nu^*_t} \ar[ddd]_{\omega^*_t}
& & & S'_t \ar[ddd]^{\omega'^*_t} \\
& \sTop(S) \ar[r]^-{\sTop(\nu)^*} \ar[lu]_{\sTop(S/\rho^t)^*}
\urtwocell\omit{}
\drtwocell\omit{_\ \ \ \ \ \ \ \ \ \sTop(\Xi)^*}
\ar[d]_{\sTop(\omega)^*} & \sTop(S') \ar[d]^{\sTop(\omega')^*}
\ar[ru]^{\sTop(S'/\rho^t)^*} \drtwocell\omit{} \\
\urtwocell\omit{}
& \sTop(T) \ar[r]_-{\sTop(\mu)^*} \ar[ld]^{\sTop(T/\rho^t)^*} &
\sTop(T') \ar[rd]_{\sTop(T'/\rho^t)^*} & \\
T_t \ar[rrr]_-{\mu^*_t} & & \ultwocell\omit{} & T'_t.
}$$
whose four unmarked $2$-cells are identities. We complete
it by adding the orientation
$$
\Xi^\dagger_t:=
(\Xi_t,\eta^{(\omega\circ\mu)_t},\eta^{(\nu\circ\omega')_t})^\dagger:
(\nu\circ\omega')^*_t\Rightarrow(\omega\circ\mu)^*_t
$$
for the front face (notation of remark
\ref{rem_adjoint-transf}(ii)). Then a direct inspection
shows that the resulting cubical diagram commutes both
on $1$-cells and $2$-cells. As in \eqref{subsec_flowers},
there follows an essentially commutative oriented diagram
of $2$-categories :
$$
\xymatrix@R+5pt{ \Stack(\Can\,T'_t)
\ar[rrrrr]^-{\St(\mu^*_t)_*} \ar[ddd]_{\St(\omega'^*_t)_*}
\ar[rd] & & & & & \Stack(\Can\,T_t)
\ar[ddd]^{\St(\omega^*_t)_*} \ar[ld] \\
\drtwocell\omit{} & \Stack(\sS(T')) \ar[rrr]^-{\St(\sTop(\mu)^*)_*}
\ar[d]|{\St(\sTop(\omega')^*)_*} &
\drtwocell\omit{_\ \ \ \ \ \ \ \ \ \St(\sTop(\Xi)^*)^\gamma_*}
& \ultwocell\omit{} & \Stack(\sS(T))
\ar[d]|{\St(\sTop(\omega)^*)_*} & \\
& \Stack(\sS(S')) \ar[rrr]_-{\St(\sTop(\nu)^*)_*}
& & & \Stack(\sS(S)) \urtwocell\omit{} & \\
\Stack(\Can\,S'_t) \ar[rrrrr]_-{\St(\nu^*_t)_*}
\ar[ru] & & \urtwocell\omit{} & & &
\Stack(\Can\,S_t) \ar[lu]
}$$
whose four unmarked $2$-cells are defined as in the diagram of
\eqref{subsec_flowers}, and which we complete to a cubical
diagram that commutes on $2$-cells, by adding the orientation
for the front face
$$
\St(\Xi^\dagger_t)^\gamma_*:\St(\omega^*_t)_*\circ\St(\mu^*_t)_*
\Rightarrow\St(\nu_t^*)_*\circ\St(\omega'^*_t)_*.
$$
Thus, let $\cE$ be any stack on the site $\Can(\sTop(T))$,
and set $\cE_t:=\St(\sTop(T/\rho^t)^*)^*\cE$; in light of
proposition \ref{prop_nerd} and remark
\ref{rem_transit-base-change}(i), we conclude that :
$$
\St(\sTop(S'/\rho^t)^*)^**\Upsilon(\St(\sTop(\Xi)^*)^\gamma_*)_\cE
\text{ is an equivalence}\Leftrightarrow
\Upsilon(\St(\Xi^\dagger_t)^\gamma_*)_{\cE_t}
\text{ is an equivalence}.
$$

\sset\subsubsection{}\label{subsec_tapeworm}
Next, consider an oriented diagram of fibred lex-sites over $I$ :
$$
\cD\quad :\quad
{\diagram \underline\cA \ar[r]^-{\phi'} \ar[d]_{\psi'}
\drtwocell\omit{_\ \beta} &
\underline\cA' \ar[d]^\psi \\
\underline\cB \ar[r]_-\phi & \underline\cB'.
\enddiagram}
$$
For every $t\in\Ob(I)$ the restriction $\cD_t$ of the diagram
$\cD$ to the fibre categories over $t$ yields another oriented
square of $2$-categories
$$
\St(\cD_t)_*\quad :\quad
{\spreaddiagramcolumns{+50pt}\diagram
\Stack(\cA_t) \ar[r]^-{\St(\phi'_t)_*} \ar[d]_{\St(\psi'_t)_*}
\drtwocell\omit{_\qquad\St(\beta_t)^\gamma_*}
& \Stack(\cA'_t) \ar[d]^{\St(\psi_t)_*} \\
\Stack(\cB_t) \ar[r]_-{\St(\phi_t)_*} & \Stack(\cB'_t).
\enddiagram}$$
For every fibred site $\underline\cA:=(\cA,p,J_\bullet)$ over $I$,
let also $i_{\cA,t}:\cA_t\to\cA$ be the inclusion functor.

\begin{corollary}\label{cor_breath-again}
In the situation of \eqref{subsec_tapeworm}, let $\cE$ be
any stack on the total site of $\underline\cA'$, and set
$\cE_t:=\St(i_{\cA',t})_*\cE$ for every $t\in\Ob(I)$.
The following conditions are equivalent :
\begin{enumerate}
\alphaenu
\item
$\Upsilon(\St(\beta)^\gamma_*)_\cE$
is an equivalence.
\item
$\Upsilon(\St(\beta_t)^\gamma_*)_{\cE_t}$
is an equivalence for every $t\in\Ob(I)$.
\end{enumerate}
\end{corollary}
\begin{proof} By invoking the pseudo-natural transformation
$\sbb_\bullet:\sS\circ\underline{\lex.\sT}\to\totSite$
of remark \ref{rem_b-bullet}(ii), we get an oriented cubical
diagram of fibred sites :
$$
\xymatrix{\sS(\underline\sT\,\underline\cA)
\ar[rrr]^-{\sS(\underline\sT\,\phi')}
\ar[ddd]_{\sS(\underline\sT\,\psi')}
\ar[rd]_{\sbb_{\underline\cA}} & & &
\sS(\underline\sT\,\underline\cA')
\ar[ddd]^{\sS(\underline\sT\,\psi)}
\ar[ld]^{\sbb_{\underline\cA'}} \\
\drtwocell\omit{^\tau_1\ } & \totSite(\underline\cA)
\ar[r]^-{\phi'} \ar[d]_{\psi'} \drtwocell\omit{_\ \beta}
& \ultwocell\omit{^\ \tau_2} \totSite(\underline\cA')
\ar[d]^\psi & \\
& \totSite(\underline\cB) \ar[r]_-\phi &
\totSite(\underline\cB') \urtwocell\omit{_\ \tau_3} \\
\sS(\underline\sT\,\underline\cB)
\ar[rrr]_-{\sS(\underline\sT\,\phi)}
\ar[ru]^{\sbb_{\underline\cB}} & \urtwocell\omit{_\ \tau_4} & &
\ar[lu]_{\sbb_{\underline\cB'}} \sS(\underline\sT\,\underline\cB')
}$$
whose front face is oriented by
$\sS(\underline\sT\beta)=\sTop(\underline\sT\beta)^*$,
and where the orientations $\tau_1,\dots,\tau_4$ are given
by the coherence constraints of $\sbb_\bullet$; again, this
cubical diagram commutes on $2$-cells. After applying
termwise as usual the pseudo-functor $\St(-)_*$, we deduce
a similar cubical diagram whose front and back faces are
oriented by $\St(\sTop(\underline\sT\beta)^*)^\gamma_*$
and respectively $\St(\beta)^\gamma_*$, and whose
other four faces are oriented by $\St(\tau_i)^\gamma_*$ for
$i=1,\dots,4$.

\begin{claim}\label{cl_muschiato}
$\St(\sbb_{\underline\cA})_*$ is a $2$-equivalence for every
fibred lex-site $\underline\cA$ over $I$.
\end{claim}
\begin{pfclaim} By construction, $\sbb_{\underline\cA}$ is the
composition of an isomorphism of canonical sites
$\sS(\underline\sT\,\underline\cA)\isom
\Can\circ\sT(\totSite\,\underline\cA)$ with the counit
$\eps_{\totSite(\underline\cA)}:\Can\circ\sT(\totSite\,\underline\cA)
\to\totSite(\underline\cA)$. Thus, it suffices to check that
$\St(\eps_{\totSite(\underline\cA)})_*$ is a $2$-equivalence.
But the natural isomorphism
$\sT(\totSite\,\underline\cA)\isom\totSite(\underline\cA)^\sim$
identifies $\eps_{\totSite(\underline\cA)}$ with the Yoneda
morphism $h^a_{\totSite(\underline\cA)}:\totSite(\underline\cA)^\sim\to
\totSite(\underline\cA)$, so the assertion follows from theorem
\ref{th_2-equiv-for-cats-of-stacks}.
\end{pfclaim}

\begin{claim}\label{cl_No-to-GMO}
For every fibred lex-site $\underline\cA$ over $I$ we have
pseudo-natural equivalences :
$$
\St(\omega^*_{\cA_t})^*\circ
\St(\sTop(\underline\sT(\underline\cA)/\rho^t)^*)^*\circ
\St(\sbb_{\underline\cA})^*\isom\St(h^a_{\cA_t})^*\circ\St(i_{\cA,t})_*
\qquad
\text{for every $t\in\Ob(I)$}
$$
where $\omega_{\cA_t}:\cA_t^\sim\isom\sT(\cA_t)$ is the
isomorphism of topoi as in \eqref{subsec_T-and-Can}.
\end{claim}
\begin{pfclaim} Recall that $\St(h^a_{\cA_t})^*$ and its $2$-adjoint
$\St(h^a_{\cA_t})_*$ are $2$-equivalences (theorem
\ref{th_2-equiv-for-cats-of-stacks}); then, from proposition
\ref{prop_breve-for-stacks}(i,ii) we get pseudo-natural
equivalences :
$$
\St(h^a_{\cA_t})^*\circ\St(i_{\cA,t})_*\isom
\St(h^a_{\cA_t})^*\circ\St(\breve\imath_{\cA,t})^*\isom
\St(\breve\imath_{\cA,t}^*)^*\circ\St(h^a_{\totSite(\underline\cA)})^*.
$$
We are thus reduced to checking the essential commutativity
of the following diagram of sites :
$$
\xymatrix@C+20pt{
\Can(\cA_t^\sim) \ar[r]^-{\breve\imath^*_{\underline\cA,t}}
\ar[d]_{\sTop(\underline\sT(\underline\cA)/\rho^t)^*\circ\omega^*_{\cA_t}} &
\Can(\totSite\,\underline\cA)^\sim
\ar[d]^{h^a_{\totSite(\underline\cA)}} \\
\Can\circ\sTop(\underline\sT\,\underline\cA)
\ar[r]^-{\sbb_{\underline\cA}} & \totSite(\underline\cA).
}$$
The latter follows by direct inspection.
\end{pfclaim}

Now, on the one hand, it is clear that condition (a) holds
if and only if $\St(i_{\cB,t})_**\Upsilon(\St(\beta)^\gamma_*)_\cE$
is an equivalence for every $t\in\Ob(I)$. By claim
\ref{cl_No-to-GMO}, the latter holds if and only if :
\begin{enumerate}
\alphaenu\addenu\addenu
\item
$(\St(\sTop(\underline\sT(\underline\cB)/\rho^t)^*)^*\circ
\St(\sbb_{\underline\cB})^*)*\Upsilon(\St(\beta)^\gamma_*)_\cE$
is an equivalence for every $t\in\Ob(I)$.
\end{enumerate}
Set $\cE':=\St(\sbb_{\underline\cA'})^*(\cE)$. By claim
\ref{cl_muschiato} and remarks \ref{rem_when-Upsilon-inverts}
and \ref{rem_transit-base-change}(i), condition (c)
in turn holds if and only if :
\begin{enumerate}
\alphaenu\addenu\addenu\addenu
\item
$\St(\sTop(\underline\sT(\underline\cB)/\rho^t)^*)^**
\Upsilon(\St(\sTop(\underline\sT\beta)^*)^\gamma_*)_{\cE'}$
is an equivalence for every $t\in\Ob(I)$.
\end{enumerate}
Then, the discussion of \eqref{subsec_first-for-topoi} shows
that (d) holds if and only if :
\begin{enumerate}
\alphaenu\addenu\addenu\addenu\addenu
\item
$\Upsilon(\St((\sT\beta_t)^\dagger
)^\gamma_*)_{\St(\sTop(\underline\sT(\underline\cA')/\rho^t)^*)^*\cE'}$ is
an equivalence for every $t\in\Ob(I)$.
\end{enumerate}
By invoking again claim \ref{cl_No-to-GMO}, we see that (e)
in turn holds if and only if :
\begin{enumerate}
\alphaenu\addenu\addenu\addenu\addenu\addenu
\item
$\Upsilon(\St((\beta^\sim_t)^\dagger)^\gamma_*)_{\St(h^a_{\cA'_t})^*\cE_t}$
is an equivalence for every $t\in\Ob(I)$.
\end{enumerate}
However, from the counit $\eps:\Can\circ\sT\Rightarrow\one_\Site$
of the $2$-adjoint pair $(\Can,\sT)$ and the system of isomorphisms
$\omega_C:C^\sim\isom\sT(C)$ of \eqref{subsec_T-and-Can} we get as
usual an oriented cubical diagram
$$
\xymatrix@C+20pt{
\Stack(\Can\,\cA^\sim_t) \ar[rrr]^-{\St(\tilde\phi'^*_t)_*}
\ar[ddd]_{\St(\tilde\psi'^*_t)_*} \ar[rd]^-{\St(h^a_{\cA_t})_*} & & &
\Stack(\Can\,\cA'^\sim_t) \ar[ddd]^{\St(\tilde\psi^*_t)_*}
\ar[ld]_-{\St(h^a_{\cA'_t})_*} \\
\drtwocell\omit{} &
\Stack(\cA_t) \ar[r]^-{\St(\phi'_t)_*} \ar[d]_{\St(\psi'_t)_*}
\drtwocell\omit{_\qquad\St(\beta_t)^\gamma_*} & \ultwocell\omit{}
\Stack(\cA'_t) \ar[d]^{\St(\psi_t)_*} & \\
& \Stack(\cB_t) \ar[r]_-{\St(\phi_t)_*} & \urtwocell\omit{}
\Stack(\cB'_t) \\
\Stack(\Can\,\cB^\sim_t) \ar[rrr]_-{\St(\tilde\phi^*_t)_*}
\ar[ru]^-{\St(h^a_{\cB_t})_*} & \urtwocell\omit{} & &
\Stack(\Can\,\cB'^\sim_t) \ar[lu]_-{\St(h^a_{\cB'_t})_*}.
}$$
commuting on $2$-cells, whose four unmarked orientations are
pseudo-natural equivalences, and whose front face is oriented
by $\St((\beta_t^\sim)^\dagger)^\gamma_*$. Combining with remarks
\ref{rem_when-Upsilon-inverts} and
\ref{rem_transit-base-change}(i) we finally conclude that
(f)$\Leftrightarrow$(b).
\end{proof}

\begin{example}\label{ex_sMorph-is-back}
Let $u:C':=(\cC',J')\to C:=(\cC,J)$ be a morphism of lex-sites.

(i)\ \
Recall that the target functor $\st_\cC:\sMorph(\cC)\to\cC$ is
a fibration whose fibre category $\st_\cC^{-1}X$ is naturally
identified with $\cC/X$, for every $X\in\Ob(\cC)$ (example
\ref{ex_fibred-cats}(iii)). Fix a unital cleavage for this
fibration, and let $\sc$ be the associated pseudo-functor;
it is easily seen that for every morphism $\phi:X\to Y$ in
$\cC$, the functor $\sc_\phi:\cC/Y\to\cC/X$ is right adjoint
to $\phi_*:\cC/X\to\cC/Y$, hence it is left exact : explicitly,
$\sc_\phi$ assigns to every object $Z\xrightarrow{h}Y$ of
$\cC/Y$ the induced projection $(X\times_YZ\to X)\in\Ob(\cC/X)$,
where $X\times_YZ$ is a choice of representative for the fibre
product of $X$ and $Z$ over $Y$ (detail left to the reader).
This description implies easily that $\sc_\phi$ is a morphism
of lex-sites $C/X\to C/Y$, for every such $\phi$ (notation of
\eqref{sec_Localization-topoi}). The collection of sites
$(C/X~|~X\in\Ob(\cC))$ then endows the fibration $\st_\cC$ with
a well defined structure of fibred lex-site
$(\sMorph(\cC),\st_\cC,J^\cC_\bullet)$.

(ii)\ \
The functor $u:\cC\to\cC'$ induces a commutative diagram of
categories :
$$
\xymatrix@C+30pt{
\sMorph(\cC) \ar[r]^-{\sMorph(u)} \ar[d]_{\st_\cC} &
\sMorph(\cC') \ar[d]^{\st_{\cC'}} \\
\cC \ar[r]^-u & \cC'.
}$$
Let now $\cB$ be any category, and $F:\cB\to\cC$ a given functor;
we set
$$
\cA:=\Fib(F)^*(\sMorph(\cC))
\qquad
\cA':=\Fib(F\circ u)^*(\sMorph(\cC'))
$$
and let
$$
\sMorph(\cC)\xleftarrow{\pi}\cA\xrightarrow{\phi}\cB
\xleftarrow{\phi'}\cA'\xrightarrow{\pi'}\sMorph(\cC')
$$
be the natural projections. There is a unique functor
$$
g:\cA\to\cA'
\qquad
\text{such that $\phi'\circ g=\phi$ and
$\pi'\circ g=\sMorph(u)\circ\pi$}
$$
and since $u$ is left exact, it is easily seen that
$g$ is cartesian (details left to the reader). Moreover,
the restriction $\phi^{-1}B\to\phi'^{-1}B$ of $g$ is naturally
identified with $u_{|FB}:\cC/FB\to\cC'/uFB$, which is a morphism
of lex-sites $C'/uFB\to C/FB$. Thus, $g$ is a morphism of
fibred lex-sites
$$
\cB\times_{\cC'}(\sMorph(\cC'),\st_{\cC'},J^{\cC'}_\bullet)
\to\cB\times_\cC(\sMorph(\cC),\st_\cC,J^\cC_\bullet)
$$
and therefore, $g$ is also a morphism of the induced total
sites (proposition \ref{prop_actually-morph-of-sites})
$$
g:(\cA',J_{\cA'})\to(\cA,J_\cA).
$$

(iii)\ \
Suppose moreover that $\cB$ is cofiltered, in which case
the pseudo-functor $\sc^o$ factors through a pseudo-functor
$\tilde\sc:\cB\to\lex.\Site$ as in
\eqref{subsec_localization-morph-of-sites}, whose strong
$2$-limit is represented by the localization $\cA[\Sigma^{-1}_\cA]$,
endowed with a certain topology $J^*_\cA$, where $\Sigma_\cA$ denotes
the set of cartesian morphisms of $\cA$. Likewise, the corresponding
localization $\cA'[\Sigma^{-1}_{\cA'}]$ of $\cA'$ carries a natural
topology $J^*_{\cA'}$, and for every $B\in\Ob(\cB)$ the composition
of the inclusion functor $i_{FB}:\cC/FB\to\cA$ (resp.
$i_{uFB}:\cC'/uFB\to\cA'$) with the localization functor
$L_\cA:\cA\to\cA[\Sigma^{-1}_\cA]$ (resp.
$L_{\cA'}:\cA'\to\cA'[\Sigma^{-1}_{\cA'}]$) is a morphism of sites
$$
l_{FB}:(\cA[\Sigma^{-1}_\cA],J^*_\cA)\to C/FB
\qquad
\text{(resp.
$l_{uFB}:(\cA'[\Sigma^{-1}_{\cA'}],J^*_{\cA'})\to C'/uFB$\ )}.
$$
Hence, for every $B\in\Ob(\cB)$ we get an oriented square
of sites :
$$
\xymatrix@C+20pt{
(\cA'[\Sigma^{-1}_{\cA'}],J^*_{\cA'}) \ar[r]^-{l_{uFB}} \ar[d]_{g_\Sigma}
\drtwocell\omit & C'/uFB \ar[d]^{u_{|FB}} \\
(\cA[\Sigma^{-1}_\cA],J^*_\cA) \ar[r]^-{l_{FB}} & C/FB
}$$
whose orientation is the identity $\one_{l_{FB}\circ g_\Sigma}$,
and where $g_\Sigma$ is the localization of $g$ (see
\eqref{subsec_bc-for-limit-site}). With this notation, we have
the following :
\end{example}

\begin{corollary}\label{cor_sMorph-is-back}
In the situation of example {\em\ref{ex_sMorph-is-back}(iii)},
suppose furthermore that $\cB$ admits a final object $B_0$,
and that every object of\/ $\cC'$ is quasi-compact for the
topology $J'$. Then the base change transformation
$$
\Upsilon(\St(\one_{l_{FB_0}\circ g_\Sigma})^\gamma_*):
\St(l_{FB_0})^*\circ\St(u_{|FB_0})_*\to
\St(g_\Sigma)_*\circ\St(l_{uFB_0})^*
$$
is a pseudo-natural equivalence.
\end{corollary}
\begin{proof} We regard $C/FB_0$ as a fibred lex-site over the
category $\bone$ with one object and one morphism; then we
may define the fibred lex-site $\cB\times_\bone C/FB_0$ over
$\cB$, as in \eqref{subsec_pullback-fib-site} : its underlying
category is $\cB\times\cC/FB_0$, and its fibre categories are
all naturally identified with $\cC/FB_0$, and endowed with the
topology of $C/FB_0$. Likewise we define the fibred lex-site
$\cB\times_\bone C'/uFB_0$, and we endow as well the underlying
categories with their respective total site topologies. According
to example \ref{ex_base-with-fin-obj}, the inclusion functors
$i_{FB_0}:\cC/FB_0\to\cB\times\cC/FB_0$ and
$i_{uFB_0}:\cC'/uFB_0\to\cB\times\cC'/uFB_0$ are then morphisms
of sites :
$$
i_{FB_0}:\cB\times_\bone C/FB_0\to C/FB_0
\qquad
i_{uFB_0}:\cB\times_\bone C'/uFB_0\to C'/uFB_0.
$$
For every $B\in\Ob(\cB)$, let $t_B:B\to B_0$ be the unique
morphism in $\cB$; we have a natural isomorphism of categories
$$
(\cC/FB_0)/Ft_B\isom\cC/FB
$$
that identifies the source functor
$\ss_{Ft_B}:(\cC/FB_0)/Ft_B\to\cC/FB_0$ with the functor
$(Ft_B)_*:\cC/FB\to\cC/FB_0$, and likewise for the category
$(\cC'/uFB_0)/uFt_B$. According to \eqref{subsec_pashimotta},
for every such $B$ we have then an oriented square of sites :
$$
\xymatrix@C+20pt{ C'/uFB \ar[r]^-{p_{uFt_B}} \ar[d]_{u_{|FB}}
\drtwocell\omit{_\ \ \ \beta_B} &
C'/uFB_0 \ar[d]^{u_{|FB_0}} \\
C/FB \ar[r]^-{p_{Ft_b}} & C/FB_0.
}$$
It is easily seen that the system of orientations
$(\beta_B~|~B\in\Ob(\cB))$ amounts to an orientation $\beta$
for the central square subdiagram in the following diagram
of oriented squares of sites :
$$
\xymatrix@C+20pt{ (\cA'[\Sigma^{-1}_{\cA'}],J^*_{\cA'}) \ar[r]^-{L_{\cA'}}
\ar[d]_{g_\Sigma} \drtwocell\omit & (\cA',J_{\cA'}) \ar[d]^g
\ar[r]^-{p_{\cA'}} \drtwocell\omit{_\ \ \beta}
& \cB\times_\bone C'/uFB_0 \ar[d]^{\cB\times u_{|FB_0}}
\ar[r]^-{i_{uFB_0}} \drtwocell\omit & C'/uFB_0 \ar[d]^{u_{|FB_0}} \\
(\cA[\Sigma^{-1}_\cA],J^*_\cA) \ar[r]^-{L_\cA} & (\cA,J_\cA)
\ar[r]^-{p_\cA} & \cB\times_\bone C/FB_0 \ar[r]^-{i_{FB_0}} & C/FB_0.
}$$
Here $p_\cA$ is the cartesian functor whose restriction to
fibre categories $\cC/FB_0\to\cC/FB$ agrees with $p_{Ft_B}$,
for every $B\in\Ob(\cB)$, and the orientations of the left
and right square are identities. Under our assumptions, it
is easily seen that every object of $\cA$ is quasi-compact
for the topology $J_{\cA'}$; then the base change transformation
associated with the left square subdiagram is a pseudo-natural
equivalence, by virtue of proposition \ref{prop_bc-lim-and-tot}.
The same holds for the base change transformation associated
with the central square, by proposition \ref{prop_pashimotta}
and corollary \ref{cor_breath-again}. A simple inspection
shows that the composition of the three top (resp. bottom)
horizontal arrows equals $l_{uFB_0}$ (resp. $l_{FB_0}$). It
remains therefore only to check that the base change transformation
$$
\St(i_{FB_0})^*\St(u_{|FB_0})_*\to\St(\cB\times u_{|FB_0})_*\St(i_{uFB_0})^*
$$
associated with the right square is a pseudo-natural equivalence.
To this aim, we remark :

\begin{claim}\label{cl_birdies}
Let $\cF_\bullet$ be any sheaf of categories on $C'/uFB_0$ such
that $\cFib(\cF_\bullet)$ is a stack on $C'/uFB_0$. Then
$\cFib((\tilde\imath{}^{\ *}_{uFB_0},\bCat)^*\cF_\bullet)$ is
a stack on $\cB\times_\bone C'/uFB_0$.
\end{claim}
\begin{pfclaim} notice first that, since $B_0$ is a final object
in $\cB$, the projection $q_{uFB_0}:\cB\times\cC'/uFB_0\to\cC'/uFB_0$
is a left adjoint for the functor $i_{uFB_0}$. On the other hand,
proposition \ref{prop_sheaves-on-tot-site}(iv) easily
implies that $q_{uFB_0}$ is continuous for the sites $C'/uFB_0$
and $\cB\times_\bone C'/uFB_0$. Combining with lemma
\ref{lem_needed}(ii), we deduce a natural isomorphism of
functors :
\set\begin{equation}\label{eq_plusieurs}
\tilde\imath{}^{\ *}_{uFB_0}\isom\tilde q_{uFB_0*}.
\end{equation}
Thus, it suffices to check that
$\cFib((\tilde q_{uFB_0*},\bCat)_*\cF_\bullet)$ is a stack. The latter
follows easily from corollary \ref{cor_lims-and-totsites}(ii).
\end{pfclaim}

In light of claim \ref{cl_birdies}, corollary
\ref{cor_conditions-for-trivial-bc}(iv), and the discussion
of \eqref{subsec_from-beta-dagger-to-tilde}, we are then
reduced to checking that the base change transformation
$$
\tilde\imath{}^{\ *}_{FB_0}\circ(u_{|FB_0})^\sim_*\to
(\cB\times u_{|FB_0})^\sim_*\circ\tilde\imath{}^{\ *}_{uFB_0}
$$
is an isomorphism of functors. But notice that $i_{uFB_0}$ is
fully faithful, and it is both continuous and cocontinuous
(proposition \ref{prop_sheaves-on-tot-site}(iii)); it follows
that $\tilde\imath{}^{\ *}_{uFB_0}$ is fully faithful (lemma
\ref{lem_improve}(iv)), hence the unit of the adjoint pair
$(\tilde\imath{}^{\ *}_{uFB_0},\tilde\imath_{uFB_0*})$ is an
isomorphism. Thus, let $\eta$ and $\eps$ be the unit and
counit for the adjoint pair
$(\tilde\imath{}^{\ *}_{FB_0},\tilde\imath_{FB_0*})$; we are
reduced to checking that
$\eps*((\cB\times u_{|FB_0})^\sim_*\circ\tilde\imath{}^{\ *}_{uFB_0})$
is an isomorphism. However, let
$q_{FB_0}:\cB\times\cC/FB_0\to\cC/FB_0$ be the projection;
arguing as in the proof of claim \ref{cl_birdies} we get
an isomorphism of functors
\set\begin{equation}\label{eq_last-plusieurs}
\tilde\imath{}^{\ *}_{FB_0}\isom\tilde q_{FB_0*}.
\end{equation}
From \eqref{eq_plusieurs} and \eqref{eq_last-plusieurs} there
follow isomorphisms of functors :
$$
(\cB\times u_{|FB_0})^\sim_*\circ\tilde\imath{}^{\ *}_{uFB_0}\isom
(\cB\times u_{|FB_0})^\sim_*\circ\tilde q_{uFB_0*}\isom
\tilde q_{FB_0*}\circ(u_{|FB_0})^\sim_*\isom
\tilde\imath{}^{\ *}_{FB_0}\circ(u_{|FB_0})^\sim_*.
$$
So it suffices to check that
$\eps*(\tilde\imath{}^{\ *}_{FB_0}\circ(u_{|FB_0})^\sim_*)$ is an
isomorphism; but again, arguing as in the foregoing we see
that $\tilde\imath{}^{\ *}_{FB_0}$ is fully faithful, so $\eta$
is an isomorphism, and finally the assertion follows, taking
into account the triangular identities of \eqref{subsec_adj-pair}.
\end{proof}

\begin{remark}\label{rem_never-ending-details}
In the situation of example \ref{ex_sMorph-is-back}(ii), the
functors $\phi:\cA\to\cB$ and $\ss_\cC\circ\pi:\cA\to\cC$
induce a functor
$$
s:\cA\to\cB\times\cC
\qquad
(B,f:X\to Y)\mapsto(B,X)
$$
(namely, $s$ is the unique functor whose composition with
the projections $\cB\leftarrow\cB\times\cC\to\cC$ equals
respectively $\phi$ and $\ss_\cC\circ\pi$). The functor $s$
admits a right adjoint :
$$
q:\cB\times\cC\to\cA
$$
that assigns to every $(B,X)\in\Ob(\cB\times\cC)$ the
object $(B,p_{X,B}:X\times FB\to FB)\in\Ob(\cA)$, where
$p_{X,B}$ is the natural projection. To every morphism
$(\beta,f):(B,X)\to(B',X')$ of $\cB\times\cC$, the functor
$q$ assigns the morphism $q(\beta,f):=(\beta,\cD)$ of $\cA$,
with $\cD$ the commutative square :
$$
\xymatrix@C+20pt{
X\times FB \ar[r]^-{p_{X,B}} \ar[d]_{f\times F\beta}
& FB \ar[d]^{F\beta} \\
X'\times FB' \ar[r]^-{p_{X',B'}} & FB'.
}$$
Let us regard $(\cC,J)$ as a fibred site over the category
$\bone$ with one object and one morphism; then
$\cB\times_\bone(\cC,J)$ is a fibred category over $\cB$, and
taking into account remark \ref{rem_continue-local}(iii), it
is easily seen that $q$ is a morphism of fibred sites
$$
q:\cB\times_\cC(\sMorph(\cC),\st_\cC,J^\cC_\bullet)\to
\cB\times_\bone(\cC,J)
$$
and therefore it is as well a morphism of the respective
total sites
$$
q:(\cA,J_\cA)\to(\cB\times\cC,J_{\cB\times\cC}).
$$
It follows that $s$ is cocontinuous for the topologies $J_\cA$
and $J_{\cB\times\cC}$ (lemma \ref{lem_needed}(i)); moreover, we
easily deduce from example \ref{ex_localize-is-weak} that $s$
is a weak morphism of fibred sites
$$
s:\cB\times_\bone(\cC,J)\to
\cB\times_\cC(\sMorph(\cC),\st_\cC,J^\cC_\bullet)
$$
(see definition \ref{def_weak-fibred-morph}); therefore it is
as well a weak morphism of the respective total sites
$$
s:(\cB\times\cC,J_{\cB\times\cC})\to(\cA,J_\cA)
$$
(proposition \ref{prop_fibred-weak-is-weak}). Combining with
proposition \ref{prop_breve-for-stacks}(ii,iii), we deduce a
pseudo-natural equivalence of pseudo-functors :
$$
\St(s)_*\isom\St(q)^*.
$$
\end{remark}

\section{Monoids and polyhedra}\label{chap_monoids}
Unless explicitly stated otherwise, {\em every monoid encountered
in this chapter shall be commutative}. For this reason, we
shall usually economize adjectives, and write just ``monoid''
when referring to commutative monoids.

\subsection{Monoids}\label{sec_toric}
If $M$ is any monoid, we shall usually denote the composition law
of $M$ by multiplicative notation: $(x,y)\mapsto x\cdot y$ (so $1$
is the neutral element). However, sometimes it is convenient to be
able to switch to an additive notation; to allow for that, we shall
denote by $(\log M,+)$ the monoid with additive composition law, whose
underlying set is the same as for the given monoid $(M,\cdot)$, and
such that the identity map is an isomorphism of monoids (then, the
neutral element of $\log M$ is denoted by $0$). For emphasis, we may
sometimes denote by $\log:M\isom\log M$ the identity map, so that
one has the tautological identities :
$$
\log 1=0\qquad\text{and}\qquad \log(x\cdot y)=\log x+\log
y\quad\text{for every $x,y\in M$}.
$$
Conversely, if $(M,+)$ is a given monoid with additive composition
law, we may switch to a multiplicative notation by writing
$(\exp M,\cdot)$, in the same way.

\sset\subsubsection{}\label{subsec_toric}
For any monoid $M$, and any two subsets $S,S'\subset M$, we let :
$$
S\cdot S':=\{s\cdot s'~|~s\in S, s'\in S'\}
$$
and $S^a$ is defined recursively for every $a\in\N$, by the rule :
$$
S^0:=\{1\}
\qquad\text{and}\qquad
S^a:=S\cdot S^{a-1}
\quad
\text{if $a>0$}.
$$
Notice that the pair $(\cP(M),\cdot)$ consisting of the set of all
subsets of $M$, together with the composition law just defined, is
itself a monoid : the neutral element is the subset $\{1\}$. In the
same vein, the exponential notation for subsets of $M$ becomes a
multiplicative notation in the monoid
$(\log\cP(M),+)=(\cP(\log M),+)$, {\em i.e.} we have the
tautological identity : $\log S^a=a\cdot\log S$, for every
$S\in\cP(M)$ and every $a\in\N$.

Furthermore, for any two monoids $M$ and $N$, the set
$\Hom_\Mnd(M,N)$ is naturally a monoid. The composition law assigns
to any two morphisms $\phi,\psi:M\to N$ their product
$\phi\cdot\psi$, given by the rule :
$\phi\cdot\psi(m):=\phi(m)\cdot\psi(m)$ for every $m\in M$.

Basic examples of monoids are the set $(\N,+)$ of natural numbers,
and the non-negative real (resp. rational) numbers $(\R_+,+)$ (resp.
$(\Q_+,+)$), with their standard addition laws.

\sset\subsubsection{}\label{subsec_present-mds}
Given a surjection $X\to Y$ of monoids, it may not be possible
to express $Y$ as a quotient of $X$ -- a problem relevant to the
construction of {\em presentations\/} for given monoids, in terms
of free monoids. For instance, consider the monoid $(\Z,\odot)$
consisting of the set $\Z$ with the composition law $\odot$ such that :
$$
x\odot y:=\left\{\begin{array}{ll}
             x+y & \text{if either $x,y\geq 0$ or $x,y\leq 0$} \\
             \max(x,y) & \text{otherwise}
          \end{array}\right.
$$
for every $x,y\in\Z$. Define a surjective map
$\phi:\N^{\oplus 2}\to(\Z,\odot)$ by the rule $(n,m)\mapsto n\odot-m$,
for every $n\in\N$. Then one verifies easily that $\Ker\,\phi=\{0\}$,
and nevertheless $\phi$ is not an isomorphism.
The right way to proceed is indicated by the following :

\begin{lemma}\label{lem_present}
Every surjective map of monoids is an effective epimorphism
(in the category $\Mnd$).
\end{lemma}
\begin{proof} (See example \ref{ex_strict-epis}(v) for
the notion of effective epimorphism.) Let $\pi:M\to N$ be a
surjection of monoids. For every monoid $X$, we have a natural
diagram of sets :
$$
\xymatrix{ \Hom_\Mnd(N,X) \ar[r]^j & \Hom_\Mnd(M,X)
\ar@<-.5ex>[r]_-{p^*_2} \ar@<.5ex>[r]^-{p^*_1} &
\Hom_\Mnd(M\times_NM,X) }
$$
where $p_1,p_2:M\times_NM\to M$ are the two natural projections,
and we have to show that the map $j$ identifies $\Hom_\Mnd(N,X)$
with the equalizer of $p^*_1$ and $p^*_2$. First of all, the
surjectivity of $\pi$ easily implies that $j$ is injective. Hence,
let $\phi:M\to X$ be any map such that $\phi\circ p_1=\phi\circ p_2$;
we have to show that $\phi$ factors through $\pi$. To this aim, it
suffices to show that the map of sets underlying $\phi$ factors as
a composition $\phi'\circ\pi$, for some map of sets $\phi':N\to X$,
since $\phi'$ will then be necessarily a morphism of monoids.
However, the forgetful functor $F:\Mnd\to\Set$ commutes with fibre
products (lemma \ref{lem_forget-me-not}(iii)), and $F(\pi)$ is an
effective epimorphism, since in the category $\Set$ all surjections
are effective epimorphisms. The assertion follows.
\end{proof}

\sset\subsubsection{}
Lemma \ref{lem_present} allows to construct presentations for an
arbitrary monoid $M$, as follows. First, we choose a surjective
map of monoids $F:=\N^{(S)}\to M$, for some set $S$.
Then we choose another set $T$ and a surjection of monoids
$\N^{(T)}\to F\times_MF$. Composing with the natural projections
$p_1,p_2:F\times_MF\to F$, we obtain a diagram :
\set\begin{equation}\label{eq_free-prent-mnd}
\xymatrix{ \N^{(T)}
\ar@<.5ex>[r]^-{q_1} \ar@<-.5ex>[r]_-{q_2} &
           \N^{(S)} \ar[r] & M
}\end{equation}
which, in view of lemma \ref{lem_present}, identifies $M$ to the
coequalizer of $q_1$ and $q_2$.

\begin{definition}\label{def_several}
Let $M$ be a monoid, $\Sigma\subset M$ a subset.
\begin{enumerate}
\item
Let $(e_\sigma~|~\sigma\in\Sigma)$ be the natural basis of the free
monoid $\N^{(\Sigma)}$. We say that $\Sigma$ is a {\em system of
generators\/} for $M$, if the map of monoids $\N^{(\Sigma)}\to M$
such that $e_\sigma\mapsto\sigma$ for every $\sigma\in\Sigma$, is a
surjection.
\item
$M$ is said to be {\em finitely generated\/} if it admits a finite
system of generators.
\item
$M$ is said to be {\em fine\/} if it is integral and finitely
generated.
\item
A {\em finite presentation\/} for $M$ is a diagram such as
\eqref{eq_free-prent-mnd} that identifies $M$ to the coequalizer of
$q_1$ and $q_2$, and such that, moreover, $S$ and $T$ are finite
sets.
\item
We say that a morphism of monoids $\phi:M\to N$ is {\em finite},
if $N$ is a finitely generated $M$-module, for the $M$-module
structure induced by $\phi$.
\end{enumerate}
\end{definition}

\begin{lemma}\label{lem_finite-pres}
{\em (i)}\ \ 
Every finitely generated monoid admits a finite presentation.
\begin{enumerate}
\addenu
\item
Let $M$ be a finitely generated monoid, and $(N_i~|~i\in I)$ a
filtered family of monoids. Then the natural map :
$$
\colim_{i\in I}\Hom_\Mnd(M,N_i)\to\Hom_\Mnd(M,\colim_{i\in I}N_i)
$$
is a bijection.
\end{enumerate}
\end{lemma}
\begin{proof}(i): Let $M$ be a finitely generated monoid, and
choose a surjection $\pi:\N^{(S)}\to M$ with $S$ a finite set. We
have seen that $M$ is the coequalizer of the two projections
$p_1,p_2:P:=\N^{(S)}\times_M\N^{(S)}\to\N^{(S)}$. For every finitely
generated submonoid $N\subset P$, let $p_{1,N},p_{2,N}:N\to\N^{(S)}$
be the restrictions of $p_1$ and $p_2$, and denote by $C_N$ the
coequalizer of $p_{1,N}$ and $p_{2,N}$. By the universal property of
$C_N$, the map $\pi$ factors through a map of monoids $\pi_N:C_N\to
M$, and since $\pi$ is surjective, the same holds for $\pi_N$. It
remains to show that $\pi_N$ is an isomorphism, for $N$ large
enough. We apply the functor $M\mapsto\Z[M]$ of
\eqref{subsec_mon-to-algs}, and we derive that $\Z[M]$ is the
coequalizer of the two maps $\Z[p_1],\Z[p_2]:\Z[P]\to\Z[S]$,
{\em i.e.} $\Z[M]\simeq\Z[S]/I$, where $I$ is the ideal generated
by $\Img(\Z[p_1]-\Z[p_2])$.
Clearly $I$ is the colimit of the filtered system of analogous
ideals $I_N$ generated by $\Img(\Z[p_{1,N}]-\Z[p_{2,N}])$, for
$N$ ranging over the filtered family $\cF$ of finitely generated
submonoids of $P$.
By noetherianity, there exists $N\in\cF$ such that $I=I_N$,
therefore $\Z[M]$ is the coequalizer of $\Z[p_{1,N}]$ and
$\Z[p_{2,N}]$. But the latter coequalizer is also the same as
$\Z[C_N]$, whence the contention.

(ii): This is a standard consequence of (i). Indeed, say that
$f_1,f_2:M\to N_i$ are two morphisms whose compositions with the
natural map $N_i\to N:=\colim_{i\in I}N_i$ agree, and pick a finite
set of generators $x_1,\dots,x_n$ for $M$. For any morphism
$\phi:i\to j$ in the filtered category $I$, denote by $g_\phi:N_i\to
N_j$ the corresponding morphism; then we may find such a morphism
$\phi$, so that $g_\phi\circ f_1(x_k)=g_\phi\circ f_2(x_k)$ for
every $k\leq n$, whence the injectivity of the map in (ii). Next,
let $f:M\to N$ be a given morphism, and pick a finite presentation
\eqref{eq_free-prent-mnd}; we deduce a morphism $g:\N^{(S)}\to N$,
and since $S$ is finite, it is clear that $g$ factors through a
morphism $g_i:\N^{(S)}\to N_i$ for some $i\in I$. Set
$g'_i:=g_i\circ q_1$ and $g''_i:=g_i\circ q_2$; by assumption, after
composing $g'_i$ and of $g''_i$ with the natural map $N_i\to N$, we
obtain the same map, so by the foregoing there exists a morphism
$\phi:i\to j$ in $I$ such that $g_\phi\circ g'_i=g_\phi\circ g''_i$.
It follows that $g_\phi\circ g_i$ factors through $M$, whence the
surjectivity of the map in (ii).
\end{proof}

\begin{definition}\label{def_radical-mon}
Let $M$ be a monoid, $I\subset M$ an ideal.
\begin{enumerate}
\item
We say that $I$ is {\em principal}, if it is cyclic, when
regarded as an $M$-module.
\item
The {\em radical\/} of $I$ is the ideal $\rad(I)$ consisting of all
$x\in M$ such that $x^n\in I$ for every sufficiently large $n\in\N$.
If $I=\rad(I)$, we also say that $I$ is a {\em radical ideal}.
\item
A {\em face\/} of $M$ is a submonoid $F\subset M$ with the following
property. If $x,y\in M$ are any two elements, and $xy\in F$, then
$x,y\in F$.
\item
Notice that the complement of a face is always an ideal. We say that
$I$ is a {\em prime ideal} of $M$, if $M\!\setminus\!I$ is a face of
$M$.
\end{enumerate}
\end{definition}

\begin{proposition}\label{prop_ideals-in-fg-mon}
Let $M$ be a finitely generated monoid, and $S$ a finitely generated
$M$-module. Then we have :
\begin{enumerate}
\item
Every submodule of $S$ is finitely generated.
\item
Especially, every ideal of $M$ is finitely generated.
\end{enumerate}
\end{proposition}
\begin{proof} Of course, (ii) is a special case of (i). To show
(i), let $S'\subset S$ be an $M$-submodule, $\Sigma\subset S'$
any system of generators. Let $\cP'(\Sigma)$ be the set of all
finite subsets of $\Sigma$, and for every $A\in\cP'(\Sigma)$, denote
by $S'_A\subset S'$ the submodule generated by $A$; clearly $S'$ is
the filtered union of the family $(S'_A~|~A\in\cP'(\Sigma))$, hence
$\Z[S']$ is the filtered union of the family of $\Z[M]$-submodules
$(\Z[S']~|~S\in\cP'(\Sigma))$. Since $\Z[M]$ is noetherian and
$\Z[S]$ is a finitely generated $\Z[M]$-module, it follows that
$\Z[S'_A]=\Z[S']$ for some finite subset $A\subset\Sigma$, whence
the contention.
\end{proof}

\sset\subsubsection{}\label{subsec_sepc-of-monoid}
Let $M$ be a monoid, $(I_\lambda~|~\lambda\in\Lambda)$ any
collection of ideals of $M$; then it is easily seen that
both $\bigcup_{\lambda\in\Lambda}I_\lambda$ and
$\bigcap_{\lambda\in\Lambda}I_\lambda$ are ideals of $M$. The
{\em spectrum\/} of $M$ is the set :
$$
\Spec\,M
$$
consisting of all prime ideals of $M$. It has a natural partial
ordering, given by inclusion of prime ideals; the minimal element of
$\Spec\,M$ is the empty ideal $\emptyset\subset M$, and the maximal
element is:
$$
\fm_M:=M\!\setminus\!M^\times.
$$
If $(\fp_\lambda~|~\lambda\in\Lambda)$ is any family of prime ideals
of $M$, then $\bigcup_{\lambda\in\Lambda}\fp_\lambda$ is a prime ideal
of $M$.

\begin{definition}
Any morphism $\phi:M\to N$ of monoids induces a natural map :
$$
\phi^*:\Spec\,N\to\Spec\,M \qquad \fp\mapsto\phi^{-1}\fp
$$
of partially ordered sets. We say that $\phi$ is {\em local},
if $\phi(\fm_M)\subset\fm_N$.
\end{definition}

\begin{corollary}\label{cor_ideals-in-fg-mon}
Let $M$ be any fine and sharp monoid. The set
$\fm_M\!\setminus\!\fm_M^2$ is finite, and is the unique minimal
system of generators of $M$.
\end{corollary}
\begin{proof} It is easily seen that any system of generators
of $M$ must contain $\Sigma:=\fm_M\!\setminus\!\fm^2_M$, hence
the latter must be a finite set. On the other hand, suppose that
there exists an element $x_0\in M$ which is not contained in the
submonoid $M'$ generated by $\Sigma$. Then we may write $x_0=x_1y_1$
for some $x_1,y_1\in\fm_M$, with $x_1\notin M'$, so $x_1$ admits
a similar decomposition. Proceeding in this way, we obtain a
sequence of elements $(x_n~|~n\in\N)$ with the property that
$Mx_n\subset Mx_{n+1}$ for every $n\in\N$. We claim that
$Mx_n\neq Mx_{n+1}$ for every $n\in\N$. Indeed, if the inequality
fails for some $n\in\N$, we may write $x_{n+1}=ax_n$ for some
$a\in\N$, and on the other hand, we have by construction
$x_n=yx_{n+1}$ for some $y\in\fm_M$; summing up, we get
$x_n=yax_n$, whence $ya=1$, since $M$ is integral, therefore
$y\in M^\times$, a contradiction.

Thus, from the given $x_0$, we have produced an infinite strictly
ascending chain of ideals of $M$, which is ruled out by virtue
of proposition \ref{prop_ideals-in-fg-mon}(ii). This means that
$x_0$ cannot exist, and the corollary follows.
\end{proof}

\begin{lemma}\label{lem_localize-sharp}
Let $f_1:M\to N_1$ and $f_2:M\to N_2$ be two local morphisms of
monoids. If $N_1$ and $N_2$ are sharp, then $N_1\amalg_MN_2$ is
sharp.
\end{lemma}
\begin{proof} Let $(a,b)\in N_1\times N_2$, and suppose there
exist $c\in M$, $a'\in N_1$, $b'\in N_2$ such that
$(a,b)=(a'f_1(c),b')$ and $(1,1)=(a',f_2(c)b')$; since $N_2$ is
sharp, we deduce $f_2(c)=b'=1$, so $b=1$. Then, since $f_2$ is
local, we get $c\in M^\times$, hence $f_1(c)=1$ and  $a=a'=1$.
One argues symmetrically in case $(1,1)=(a'f_1(c),b')$ and
$(a,b)=(a',f_2(c)b')$. We conclude that $(a,b)$ represents the
unit class in $N_1\amalg_MN_2$ if and only if $a=b=1$. Now,
suppose that the class of $(a,b)$ is invertible in
$N_1\amalg_MN_2$; it follows that there exists $(c,d)$
such that $ac=1$ and $bd=1$, which implies that $a=1$ and
$b=1$, whence the contention.
\end{proof}

\begin{lemma}\label{lem_std-verif}
Let $S\subset M$ be any submonoid. The localization $j:M\to S^{-1}M$
induces an injective map $j^*:\Spec\,S^{-1}M\to\Spec\,M$ which
identifies $\Spec\,S^{-1}M$ with the subset of $\Spec\,M$ consisting
of all prime ideals $\fp$ such that $\fp\cap S=\emptyset$.
\end{lemma}
\begin{proof} For every $\fp\in\Spec\,M$, denote by $S^{-1}\fp$ the
ideal of $S^{-1}M$ generated by the image of $\fp$. We claim that
$\fp=j^*(S^{-1}\fp)$ for every $\fp\in\Spec\,M$ such that $\fp\cap
S=\emptyset$. Indeed, clearly $\fp\subset j^*(S^{-1}\fp)$; next, if
$f\in j^*(S^{-1}\fp)$, there exist $s\in S$ and $g\in\fp$ such that
$s^{-1}g=f$ in $S^{-1}M$; therefore there exists $t\in S$ such that
$tg=tsf$ in $M$, especially $tsf\in\fp$, hence $f\in\fp$, since
$t,s\notin\fp$. Likewise, one checks easily that $S^{-1}\fp$ is
a prime ideal if $\fp\cap S=\emptyset$, and $\fq=S^{-1}(j^*\fq)$
for every $\fq\in\Spec\,S^{-1}M$, whence the contention.
\end{proof}

\begin{remark}\label{rem_localize-monoids}
(i)\ \
If we take $S_\fp:=M\!\setminus\!\fp$, the complement of a prime
ideal $\fp$ of $M$, we obtain the monoid
$$
M_\fp:=S_\fp^{-1}M
$$
and $\Spec\,M_\fp\subset\Spec\,M$ is the subset consisting of all
prime ideals $\fq$ contained in $\fp$.

(ii)\ \
Likewise, if $\fp\subset M$ is any prime ideal, the spectrum
$\Spec(M\!\setminus\!\fp)$ is naturally identified with the subset
of $\Spec\,M$ consisting of all prime ideals $\fq$ containing $\fp$
(details left to the reader).

(iii)\ \
Let $S\subset M$ be any submonoid. Then there exists a smallest
face $F$ of $M$ containing $S$ (namely, the intersection of all
the faces that contain $S$). It is easily seen that $F$ is the
subset of all $x\in M$ such that $xM\cap S\neq\emptyset$. From
this characterization, it is clear that $S^{-1}M=F^{-1}M$.
In other words, every localization of $M$ is of the type $M_\fp$
for some $\fp\in\Spec\,M$.
\end{remark}

\begin{lemma}\label{lem_radical}
Let $M$ be a monoid, and $I\subset M$ any ideal. Then $\rad(I)$
is the intersection of all the prime ideals of $M$ containing $I$.
\end{lemma}
\begin{proof} It is easily seen that a prime ideal containing
$I$ also contains $\rad(I)$. Conversely, say that $f\in
M\setminus\rad(I)$; let $\phi:M\to M_f$ be the localization map.
Denote by $\fm$ the maximal ideal of $M_f$. We claim that
$I\subset\phi^{-1}\fm$, Indeed, otherwise there exist $g\in I$,
$h\in M$ and $n\in\N$ such that $g^{-1}=f^{-n}h$ in $M_f$; this
means that there exists $m\in\N$ such that $f^{m+n}=f^mgh$ in $M$,
hence $f^{m+n}\in I$, which contradicts the assumption on $f$. On
the other hand, obviously $f\notin\phi^{-1}\fm$.
\end{proof}

\begin{lemma}\label{lem_spec-quots}
{\em (i)}\ \ Let $M$ be a monoid, and $G\subset M^\times$ a subgroup.
\begin{list}{}{}
\item
\begin{enumerate}
\alphaenu
\item
The map given by the rule : $I\mapsto I/G$ establishes a natural
bijection from the set of ideals of $M$ onto the set of ideals of
$M/G$.
\item
Especially, the natural projection $\pi:M\to M/G$ induces a
bijection :
$$
\pi^*:\Spec\,M/G\to\Spec\,M
$$
\end{enumerate}
\end{list}
\begin{enumerate}
\addenu
\item
Let $(M_i~|~i\in I)$ be any finite family of monoids, and for each
$j\in I$, denote by $\pi_j:\prod_{i\in I}M_i\to M_j$ the natural
projection. The induced map
$$
\prod_{i\in I}\Spec\,M_i\to\Spec\prod_{i\in I}M_i
\quad : \quad
(\fp_i~|~i\in I)\mapsto\bigcup_{i\in I}\pi_i^*\fp_i
$$
is a bijection.
\item
Let $(M_i~|~i\in I)$ be any filtered system of monoids. The
natural map
$$
\Spec\,\colim_{i\in I}M_i\to\lim_{i\in I}\Spec\,M_i
$$
is a bijection.
\end{enumerate}
\end{lemma}
\begin{proof}(i): By lemma \ref{lem_special-p-out}(iii), $M/G$ is the
set-theoretic quotient of $M$ by the translation action of $G$. By
definition, any  ideal $I$ of $M$ is stable under the $G$-action,
hence the quotient $I/G$ is well defined, and one checks easily that
it is an ideal of $M/G$. Moreover, if $\fp\subset M$ is a prime
ideal, it is easily seen that $\fp/G$ is a prime ideal of $M/G$.
Assertions (a) and (b) are straightforward consequences.

(ii): The assertion can be rephrased by saying that every face
$F$ of $\prod_{i\in I}M_i$ is a product of faces $F_i\subset M_i$.
However, if $\underline m:=(m_i~|~i\in I)\in F$, then, for each
$i\in I$ we can write
$\underline m=\underline m(i)\cdot\underline n(i)$, where, for each
$j\in I$,  the $j$-th-component of $\underline m(i)$ (resp. of
$\underline n(i)$) equals $1$ (resp. $m_j$), unless $j=i$, in which
case it equals $m_i$ (resp. $1$). Thus, $\underline m(i)\in F$ for
every $i\in I$, and the contention follows easily.

(iii): Denote by $M$ the colimit of the system $(M_i~|~i\in I)$,
and $\phi_i:M_i\to M$ the natural morphisms of monoids, as well
as $\phi_f:M_i\to M_j$ the transition maps, for every morphism
$f:i\to j$ in $I$. Recall that the set underlying $M$ is the
colimit of the system of sets $(M_i~|~i\in I)$ (lemma
\ref{lem_forget-me-not}(iii)). Let now
$\fp_\bullet:=(\fp_i~|~\in I)$ be a compatible system of prime
ideals, {\em i.e.} such that $\fp_i\in\Spec\,M_i$ for every
$i\in I$, and $\phi_f^{-1}\fp_j=\fp_i$ for every $f:i\to j$.
We let $\beta(\fp_\bullet):=\bigcup_{i\in I}\phi_i(\fp_i)$.
We claim that $\beta(\fp_\bullet)$ is a prime ideal of
$M$. Indeed, suppose that $x,y\in M$ and $xy\in\beta(\fp_\bullet)$;
since $I$ is filtered, we may find $i\in I$, $x_i,y_i\in M_i$,
and $z_i\in\fp_i$, such that $x=\phi_i(x_i)$, $y=\phi_i(y_i)$,
and $xy=\phi_i(z_i)$. Especially, $\phi_i(z_i)=\phi_i(x_iy_i)$,
so there exists a morphism $f:i\to j$ such that
$\phi_f(z_i)=\phi_f(x_iy_i)$. But $\phi_f(z_i)\in\fp_j$, so
either $\phi_f(x_i)\in\fp_j$ or $\phi_f(y_i)\in\fp_j$, and
finally either $x\in\fp$ or $y\in\fp$, as required.

Let $\fp\subset M$ be any prime ideal; it is easily seen that
$\beta(\phi_i^{-1}\fp~|~i\in I)=\fp$. To conclude, it suffices
to show that $\fp_i=\phi_i^{-1}\beta(\fp_\bullet)$, for every
compatible system $\fp_\bullet$ as above, and every $i\in I$.
Hence, fix $i\in I$ and pick $x_i\in M_i$ such that
$\phi_i(x_i)\in\beta(\fp_\bullet)$; then there exist $j\in I$
and $x_j\in\fp_j$ such that $\phi_i(x_i)=\phi_j(x_j)$. Since
$I$ is filtered, we may find $k\in I$ and morphisms
$f:i\to k$ and $g:j\to k$ such that $\phi_f(x_i)=\phi_g(x_j)$,
so $\phi_f(x_i)\in\fp_k$, and finally $x_i\in\fp_i$, as sought.
\end{proof}

\begin{remark} In case $(M_i~|~i\in I)$ is an infinite family
of monoids, the natural map of lemma \ref{lem_spec-quots}(ii)
is still injective, but it is not necessarily surjective.
For instance, let $I$ be any infinite set, and let $\cU\subset\cP(I)$
be a non-principal ultrafilter; denote by ${}^*\N$ the quotient
of $\N^I$ under the equivalence relation $\sim_\cU$ such that
$(a_i~|~i\in I)\sim_\cU(b_i~|~i\in I)$ if and only if there
exists $U\in\cU$ such that $a_i=b_i$ for every $i\in U$.
It is clear that the composition law on $\N^I$ descends to
${}^*\N$; the resulting structure $({}^*\N,+)$ is called the
monoid of hypernatural numbers. Denote by $\pi:\N^I\to{}^*\N$
the projection, and let $\underline 0\in\N^I$ be the unit;
then it is easily seen that $\{\pi(0)\}$ is a face of ${}^*\N$,
but $\pi^{-1}(\pi(0))$ is not a product of faces $F_i\subset\N$.
\end{remark}

\begin{definition}\label{def_dim-mon}
Let $M$ be a monoid.
\begin{enumerate}
\item
The {\em dimension\/} of $M$, denoted $\dim M\in\N\cup\{+\infty\}$,
is defined as the supremum of all $r\in\N$ such that there exists a
chain of strict inclusions of prime ideals of $M$ :
$$
\fp_0\subset\fp_1\subset\cdots\subset\fp_r.
$$
\item
The {\em height\/} of a prime ideal $\fp\in\Spec\,M$ is defined as
$\hgt\,\fp:=\dim M_\fp$.
\item
A {\em facet\/} of $M$ is the complement of a prime ideal of $M$ of
height one.
\end{enumerate}
\end{definition}

\begin{remark}
(i)\ \ 
Notice that not all epimorphisms in $\Mnd$ are surjections
on the underlying sets; for instance, every localization
map $M\to S^{-1}M$ is an epimorphism.

(ii)\ \ 
If $\phi:N\to M$ is a map of fine monoids (see definition
\ref{def_several}(vi)), then it will follow from corollary
\ref{cor_fibres-are-fg} that $N\times_MN$ is also finitely
generated. If $M$ is not integral, then this fails in general : a
counter-example is provided by the morphism $\phi$ constructed in
\eqref{subsec_present-mds}.

(iii)\ \ 
Let $\Sigma$ be any set; it is easily seen that a free
monoid $M\simeq\N^{(\Sigma)}$ admits a unique minimal system of
generators, in natural bijection with $\Sigma$. Especially, the
cardinality of $\Sigma$ is determined by the isomorphism class of
$M$; this invariant is called the {\em rank} of the free monoid $M$.
This is the same as the rank of $M$ as an $\N$-module (see
example \ref{ex_rank-mod}).

(iv)\ \ 
A submonoid of a finitely generated monoid is not necessarily
finitely generated. For instance, consider the submonoid
$M\subset\N^{\oplus 2}$, with $M:=\{(0,0)\}\cup\{(a,b)~|~a>0\}$.
However, the following result shows that a face of a finitely
generated monoid is again finitely generated.
\end{remark}

\begin{lemma}\label{lem_face}
Let $f:M\to N$ be a map of monoids, $F\subset N$ a face of $N$, and
$\Sigma\subset N$ a system of generators for $N$. Then :
\begin{enumerate}
\item
$N^\times$ is a face of $N$, and $f^{-1}F$ is a face of $M$.
\item
$\Sigma\cap F$ is a system of generators for $F$.
\item
If $N$ is finitely generated, $\Spec\,N$ is a finite set, and $\dim
N$ is finite.
\item
If $N$ is finitely generated (resp. fine) then the same holds for
$F^{-1}N$.
\end{enumerate}
\end{lemma}
\begin{proof} (i) and (ii) are left to the reader, and (iii)
is an immediate consequence of (ii). To show (iv), notice that -- in
view of (ii) -- the set $\Sigma\cup\{f^{-1}~|~f\in F\cap\Sigma\}$ is
a system of generators of $F^{-1}N$.
\end{proof}

\begin{definition}
(i)\ \
If $M$ is a pointed or not pointed monoid (see remark
\ref{rem_apparent-reasons}(ii)), and $S$ is a pointed
$M$-module, we say that $S$ is {\em integral}, if for
every $x,y\in M$ and every $s\in S$ such that $xs=ys\neq 0$,
we have $x=y$. The {\em annihilator ideal\/} of $S$ is the ideal
$$
\Ann_M(S):=\{m\in M~|~ms=0\ \text{for every $s\in S$}\}.
$$
If $s\in S$ is any element, we also write $\Ass_M(s):=\Ass_M(Ms)$.
The {\em support\/} of $S$ is the subset :
$$
\Supp\,S:=\{\fp\in\Spec\,M~|~S_\fp\neq 0\}.
$$

(ii)\ \
A pointed monoid $(M,0_M)$ is called {\em integral},
if it is integral when regarded as a pointed $M$-module; it is
called {\em fine}, if it is finitely generated and integral in
the above sense.

(iii)\ \
The forgetful functor $\Mnd_\circ\to\Set$ (notation of
\eqref{subsec_wish-now}) admits a left adjoint, which assigns
to any set $\Sigma$ the {\em free pointed monoid\/}
$\N^{(\Sigma)}_\circ:=(\N^{(\Sigma)})_\circ$.

(iv)\ \
 A morphism $\phi:M\to N$ of pointed monoids is {\em local},
if both $M,N\neq 0$, and $\phi$ is local when regarded as a
morphism of non-pointed monoids.
\end{definition}

\begin{remark}\label{rem_support-mods-mon}
(i)\ \ 
Quite generally, a (non-pointed) monoid $M$ is finitely generated
(resp. free, resp. integral, resp. fine) if and only if $M_\circ$
has the corresponding property for pointed monoids. However, there
exist integral pointed monoids which are not of the form $M_\circ$
for any non-pointed monoid $M$.

(ii)\ \ 
Let $(M,0_M)$ be a pointed monoid. An {\em ideal\/} of
$(M,0_M)$ is a pointed submodule $I\subset M$. Just as for
non-pointed monoids, we say that $I$ is a {\em prime ideal}, if
$M\!\setminus\!I$ is a (non-pointed) submonoid of $M$, and a
non-pointed submonoid which is the complement of a prime ideal,
is called a {\em face\/} of $M$. Hence the smallest ideal is
$\{0\}$. Notice though, that $\{0\}$ is not necessarily a prime
ideal, hence the spectrum $\Spec\,(M,0_M)$ does not always admit
a least element. However, if $M=N_\circ$ for some non-pointed
monoid $N$, the natural morphism of monoids $N\to N_\circ$
induces a bijection :
$$
\Spec\,(N_\circ,0_{N_\circ})\to\Spec\,N.
$$

(iii)\ \
Let $I\subset M$ be any ideal; the inclusion map $I\to M$ can be
regarded as a morphism of pointed $M$-modules (if $M$ is not
pointed, this is achieved via the faithful embedding
\eqref{eq_faithf-point}), whence a pointed $M$-module $M/I$, with a
natural morphism $M\to M/I$. The latter map is also a morphism of
monoids, for the obvious monoid structure on $M/I$. One checks
easily that if $M$ is integral, $M/I$ is an integral pointed
monoid.

(iv)\ \
Let $M$ be a (pointed or not pointed) monoid, $\fp\subset M$
a prime ideal. Then the natural morphism of monoids $M\to M/\fp$
induces a bijection :
$$
\Spec\,M/\fp\isom\{\fq\in\Spec\,M~|~\fp\subset\fq\}=
\Spec\,M\!\setminus\!\fp.
$$

(v)\ \
Let $M$ be a (pointed or not pointed) monoid, and $S\neq 0$ a
pointed $M$-module. Then the support of $S$ contains at least
the maximal ideal of $M$. This trivial observation shows that
a pointed $M$-module is $0$ if and only if its support is empty.

(vi)\ \
Let $M$ be a pointed monoid, and $\Sigma\subset M$ a
{\em non-pointed\/} submonoid. The localization $\Sigma^{-1}M$
(defined in the category of monoids, as in
\eqref{subsec_localize-mons}) is actually a pointed monoid :
its zero element $0_{\Sigma^{-1}M}$ is the image of $0_M$.

(vii)\ \
If $M$ is a (pointed or not pointed) monoid, $\Sigma\subset M$
any {\em non-pointed\/} submonoid, and $S$ a pointed $M$-module,
we let as usual $\Sigma^{-1}S:=\Sigma^{-1}M\otimes_MS$ (see
remark \ref{rem_point-unpointed-topos}(i), if $M$ is not pointed).
The resulting functor $M\Mod_\circ\to \Sigma^{-1}M\Mod_\circ$
is exact. Indeed, it is right exact, since it is left adjoint
to the restriction of scalars arising from the localization
map $M\to\Sigma^{-1}M$ (see \eqref{subsec_tens-restr}), and
one verifies directly that it commutes with finite limits.
Also, it is clear that
$$
\Supp\,\Sigma^{-1}S=\Supp\,S\cap\Spec\,\Sigma^{-1}M.
$$

(viii)\ \
Let $M$ be a pointed monoid, and $N\subset M$ a pointed
submonoid. Since the final object $1$ of the category of
pointed monoids is not isomorphic to the initial object
$1_\circ$, the push-out of the diagram $1\leftarrow N\to M$
is not an interesting object (it is always isomorphic to
$1$). Even if we form the quotient $M/N$ in the category
of non-pointed monoids, we still get always $1$, since
$0_M\in N$, and therefore in the quotient $M/N$ the images
of $0_M$ and of the unit of $M$ coincide.

The only case that may give rise to a non-trivial quotient,
is when $N$ is {\em non-pointed}; in this situation we may
form $M/N$ {\em in the category of non-pointed monoids},
and then remark that the image of $0_M$ yields a zero
element $0_{M/N}$ for $M/N$, so the latter is a pointed
monoid.
\end{remark}

\begin{example}\label{ex_flat-quots}
(i)\ \
Let $M$ be a (pointed or not pointed) monoid, $G\subset M^\times$
a subgroup, and $S$ a pointed $M$-module. Then $M/G\otimes_MS=S/G$
is the set of orbits of $S$ under the induced $G$-action.

(ii)\ \
In the situation of (i), notice that the functor
$$
M\Mod\to M/G\Mod
\quad : \quad
S\mapsto S/G
$$
is exact, hence $M/G$ is a flat $M$-module. (See definition
\ref{def_pointed-flat}(i).)

(iii)\ \
Likewise, if $\Sigma\subset M$ is a non-pointed submonoid,
then the localization $\Sigma^{-1}M$ is a flat $M$-module,
due to remark \ref{rem_support-mods-mon}(vii).

(iv)\ \
Suppose that $S$ is an integral pointed $M$-module (with $M$
either pointed or not pointed), and let $\Sigma\subset M$
be a non-pointed submonoid. Then $\Sigma^{-1}S$ is also an
integral pointed $\Sigma^{-1}M$-module. Indeed, suppose that
the identity
\set\begin{equation}\label{eq_id-in-loc}
(s^{-1}a)\cdot(s'{}^{-1}b)=(s^{-1}a)\cdot(s''{}^{-1}c)\neq 0
\end{equation}
holds for some $a,b,c\in S$ and $s,s',s''\in\Sigma$; we need to
check that $s'{}^{-1}b=s''{}^{-1}c$, or equivalently, that
$s''b=s'c$ in $\Sigma^{-1}S$. However, \eqref{eq_id-in-loc}
is equivalent to $s''ab=s'ac$ in $\Sigma^{-1}S$, and the
latter holds if and only if there exists $t\in\Sigma$ such
that $ts''ab=ts'ac$ in $S$. The two sides in the latter
identity are $\neq 0$, as the same holds for the two sides
of the identity \eqref{eq_id-in-loc}; therefore $s''b=s'c$
holds already in $S$, and the contention follows.

(v)\ \
Likewise, in the situation of (iv),
$S/\Sigma:=S\otimes_MM/\Sigma$ is an integral pointed
$M/\Sigma$-module. Indeed, notice the natural identification
$S/\Sigma=
\Sigma^{-1}S\otimes_{\Sigma^{-1}M}(\Sigma^{-1}M)/\Sigma^\gp$
which -- in view of (iv) -- reduces the proof to the case
where $\Sigma$ is a subgroup of $M^\times$. Then the assertion
is easily verified, taking into account (i).

Especially, if $M$ is an integral pointed monoid, and
$\Sigma\subset M$ is any non-pointed submonoid, then
both $\Sigma^{-1}M$ and $M/\Sigma$ are integral pointed
monoids (this generalizes lemma \ref{lem_integral-quot}).

(vi)\ \
Let $G$ be any abelian group, $\phi:M\to G$ a morphism of
non-pointed monoids. Then $G_\circ$ is a flat $M_\circ$-module.
For the proof, we may -- in light of (iii) -- replace $M$
by $M^\gp$, thereby reducing to the case where $M$ is a group.
Next, by (ii), we may assume that $\phi$ is injective, in which
case $G$ is a free $M$-module with basis $G/M$.
\end{example}

\begin{remark}\label{rem_tens-is-pushout}
(i)\ \
Let $M\to N$ and $M\to N'$ be morphisms of pointed monoids; $N$
and $N'$ can be regarded as pointed $M$-modules in an obvious way,
hence we may form the tensor product $N'':=N\otimes_MN'$; the
latter is endowed with a unique monoid structure such that the
maps $e:N\to N''$ and $e':N'\to N''$ given by the rule
$n\mapsto n\otimes 1$ for all $n\in N$ (resp.
$n'\mapsto 1\otimes n'$ for all $n'\in N'$) are morphisms of
monoids. Just as for usual ring homomorphisms, the monoid $N''$
is a coproduct of $N$ and $N'$ over $M$, {\em i.e.} there is a
unique isomorphism of pointed monoids:
\set\begin{equation}\label{eq_tens-is-pushout}
N\otimes_MN'\isom N\amalg_MN'
\end{equation}
that identifies $e$ and $e'$ to the natural morphisms $N\to
N\amalg_MN'$ and $N'\to N\amalg_MN'$. As usual all this extends to
non-pointed monoids. (Details left to the reader.)

(ii)\ \
Especially, if we take $M=\{1\}_\circ$, the initial object
in $\Mnd_\circ$, we obtain an explicit description of the
pointed monoid $N\oplus N'$ : it is the quotient
$(N\times N')\!/\!\!\sim$, where $\sim$ denotes the minimal
equivalence relation such that $(x,0)\sim(0,x')$ for every
$x\in N$, $x'\in N'$. From this, a direct calculation shows
that a direct sum of pointed integral monoids is again a
pointed integral monoid.
\end{remark}

\begin{remark}\label{rem_integr-modules}
(i)\ \  
Clearly every pointed $M$-module $S$ is the colimit of the
filtered family of its finitely generated submodules. Moreover,
$S$ is the colimit of a filtered family of finitely presented
pointed $M$-modules. Recall the standard argument : pick a
countable set $I$, and let $\cC$ be the (small) full subcategory
of the category $M\Mod_\circ$ whose objects are the coequalizers
of every pair of maps of pointed $M$-modules
$p,q:M^{(I_1)}\to M^{(I_2)}$, for every finite sets
$I_1,I_2\subset I$ (this means that, for every such pair
$p,q$ we pick one representative for this coequalizer).
Then there is a natural isomorphism of pointed $M$-modules :
$$
\colim_{i\cC\!/\!S}\,\iota_S\isom S
$$
where $i:\cC\to M\Mod_\circ$ is the inclusion functor, and
$\iota_S$ is the functor as in \eqref{subsec_fibreovercat}.

(ii)\ \ 
If $S$ is finitely generated, we may find a finite filtration
of $S$ by submodules $0=S_0\subset S_1\subset\cdots\subset S_n=S$
such that $S_{i+1}/S_i$ is a cyclic $M$-module, for every
$i=1,\dots,n$.

(iii)\ \ 
Notice that if $S$ is integral and $S'\subset S$ is any
submodule, then $S/S'$ is again integral. Moreover, if $S$ is
integral and cyclic, we have a natural isomorphism of $M$-modules :
$$
S\isom M/\Ann_M(S).
$$
(Details left to the reader.)

(iv)\ \ 
Suppose furthermore, that $M^\sharp$ is finitely generated,
and $S$ is any pointed $M$-module. Lemma \ref{lem_spec-quots}(i.a)
and proposition \ref{prop_ideals-in-fg-mon}(ii) easily imply that
every ascending chain
$$
I_0\subset I_1\subset I_2\subset\cdots
$$
of ideals of $M$ is stationary; especially, the set
$\{\Ann_M(s)~|~s\in S\!\setminus\!\{0\}\}$ admits maximal elements.
Let $I$ be a maximal element in this set; a standard argument
as in commutative algebra shows that $I$ is a prime ideal :
indeed, say that $xy\in I=\Ann_M(s)$ and $x\notin I$; then
$xs\neq 0$, hence $y\in\Ann_M(xs)=I$, by the maximality of $I$.
Now, if $S$ is also finitely generated, it follows that we may
find a finite filtration of $S$ as in (ii) such that, additionally,
each quotient $S_{i+1}/S_i$ is of the form $M/\fp$, for some
prime ideal $\fp\subset M$.

These properties make the class of integral pointed modules
especially well behaved : essentially, the full subcategory
$M\text{-}\mathbf{Int.Mod}_\circ$ of $M\Mod_o$ consisting
of these modules mimics closely a category of $A$-modules
for a ring $A$, familiar from standard linear algebra. This
shall be amply demonstrated henceforth. For instance, we point
out the following combinatorial version of Nakayama's lemma :
\end{remark}

\begin{proposition}\label{prop_Naka-combi}
Let $M$ be a (pointed or not pointed) monoid, $S$ a finitely
generated integral pointed $M$-module, and $S'\subset S$ a
pointed submodule, such that
$$
S=S'\cup\fm_M\cdot S.
$$
Then $S=S'$.
\end{proposition}
\begin{proof} After replacing $S$ by $S/S'$, we may assume that
$\fm_MS=S$, in which case we have to check that $S=0$. Suppose
then, that $S\neq 0$; from remark \ref{rem_integr-modules}(ii,iii)
it follows that $S$ admits a (pointed) submodule $T\subset S$
such that $S/T\simeq M/\fm_M$. Especially, $\fm_M\cdot(S/T)=0$,
{\em i.e.} $\fm_MS\subset T$, and therefore $S=T$, which
contradicts the choice of $T$. The contention follows.
\end{proof}

\begin{remark}
(i)\ \ 
The integrality assumption cannot be omitted in proposition
\ref{prop_Naka-combi}. Indeed, take $M:=\N$ and $S:=0_\circ$,
where $0$ denotes the final $\N$-module. Then $S\neq 0$, but
$\fm_MS=S$.

(ii)\ \
Let us say that an element of the $M$-module $S$ is
{\em primitive}, if it does not lie in $\fm_MS$. We
deduce from proposition \ref{prop_Naka-combi}, the
following :
\end{remark}

\begin{corollary}
Let $M$ be a sharp (pointed or not pointed)
monoid, $S$ a finitely generated integral pointed monoid.
Then the set $S\setminus\fm_MS$ of primitive elements of $S$
is finite, and is the unique minimal system of generators
of $S$.
\end{corollary}
\begin{proof} Indeed, it is easily seen that every system
of generators of $S$ must contain all the primitive elements,
so $S\setminus\fm_MS$ must be finite. On the other hand, let
$S'\subset S$ be the submodule generated by the primitive
elements; clearly $S'\cup\fm_MS=S$, hence $S'=S$, by proposition
\ref{prop_Naka-combi}.
\end{proof}

\sset\subsubsection{}\label{subsec_LaRa}
Let $R$ be any ring, $M$ a non-pointed monoid. Notice that the
$M$-module underlying any $R[M]$-module is naturally pointed,
whence a forgetful functor $R[M]\Mod\to M\Mod_\circ$.
The latter admits a left adjoint
$$
M\Mod_\circ\to R[M]\Mod
\quad : \quad
(S,0_S)\mapsto R\La S\Ra:=\Coker\,R[0_S].
$$
Likewise, the monoid $(A,\cdot)$ underlying any (commutative unital)
$R$-algebra $A$ is naturally pointed, whence a forgetful functor
$R\Alg\to\Mnd_\circ$, which again admits a left adjoint
$$
\Mnd_\circ\to R\Alg
\quad : \quad
(M,0_M)\mapsto R\La M\Ra:=R[M]/(0_M)
$$
where $(0_M)\subset R[M]$ denotes the ideal generated by the
image of $0_M$.

If $(M,0_M)$ is a pointed monoid, and $S$ is a pointed
$(M,0_M)$-module, then notice that $R\La S\Ra$ is actually
a $R\La M\Ra$-module, so we have as well a natural functor
$$
(M,0_M)\Mod_\circ\to R\La M\Ra\Mod
\qquad
S\mapsto R\La S\Ra
$$
which is again left adjoint to the forgetful functor.

For instance, let $I\subset M$ be an ideal; from the foregoing, it
follows that $R\La M/I\Ra$ is naturally an $R[M]$-algebra, and we
have a natural isomorphism :
$$
R\La M/I\Ra\isom R[M]/IR[M].
$$
Explicitly, for any $x\in F:=M\!\setminus\!I$, let $\bar x\in R\La
M/I\Ra$ be the image of $x$; then $R\La M/I\Ra$ is a free
$R$-module, with basis $(\bar x~|~x\in F)$. The multiplication law
of $R\La M/I\Ra$ is determined as follows. Given $x,y\in F$, then
$\bar x\cdot\bar y=\bar{xy}$ if $xy\in F$, and otherwise it equals
zero.

Notice that if $I_1$ and $I_2$ are two ideals of $M$, we have a
natural identification :
$$
R\La M/(I_1\cap I_2)\Ra\isom R\La M/I_1\Ra\times_{R\La M/(I_1\cup
I_2)\Ra}R\La M/I_2\Ra.
$$
These algebras will play an important role in section
\ref{sec_log-regular}. As a special case, suppose that $\fp\subset
M$ is a prime ideal; then the inclusion $M\!\setminus\!\fp\subset M$
induces an isomorphism of $R$-algebras :
$$
R[M\!\setminus\!\fp]\isom R\La M/\fp\Ra.
$$
Furthermore, general nonsense yields a natural isomorphism
of $R$-modules :
\set\begin{equation}\label{eq_yet-more-nonsense}
R\La S\otimes_MS'\Ra\isom R\La S\Ra\otimes_{R\La M\Ra}R\La S'\Ra
\qquad
\text{for all pointed $M$-modules $S$ and $S'$.}
\end{equation}

\sset\subsubsection{}\label{subsec_restrict-scalar-mon}
Let $A$ be a commutative ring with unit, and $f:M\to(A,\cdot)$
a morphism of pointed monoids. Then $f$ induces (forgetful) functors :
$$
A\Mod\to M\Mod_\circ \qquad A\Alg\to M/\Mnd_\circ
$$
(notation of \eqref{subsec_slice-cat}) which admit left adjoints :
$$
\begin{aligned}
M\Mod_\circ & \to A\Mod & \qquad & S\mapsto S\otimes_MA:=\Z\La
S\Ra\otimes_{\Z\La M\Ra}A \\
M/\Mnd_\circ & \to A\Alg & \qquad & N\mapsto N\otimes_MA:=\Z\La
N\Ra\otimes_{\Z\La M\Ra}A.
\end{aligned}
$$
Sometimes we may also use the notation :
$$
S\otimes_MN:=\Z\La S\Ra\otimes_{\Z\La M\Ra}N \quad\text{and}\quad
S\derotimes_MK_\bullet:=\Z\La S\Ra\derotimes_{\Z\La M\Ra}K_\bullet
$$
for any pointed $M$-module $S$, any $A$-module $N$ and any object
$K_\bullet$ of $\sD^-(A\Mod)$. The latter derived tensor product
is obtained by tensoring with a flat $\Z\La M\Ra$-flat resolution
of $\Z\La S\Ra$. (Such resolutions can be constructed combinatorially,
starting from a simplicial resolution of $S$.) All the verifications
are standard, and shall be left to the reader.

Moreover, for any ideal $I\subset M$ and any $A$-module $N$,
we shall write $f(I)N$, or sometimes just $IN$, for the
$A$-submodule of $N$ generated by the system
$(f(x)\cdot y~|~x\in I,\ y\in N)$.

\begin{definition}\label{def_M-flatness}
Let $M$ be a pointed monoid, $A$ a commutative ring with unit,
$\phi:M\to(A,\cdot)$ a morphism of pointed monoids, $N$ an $A$-module.
We say that $N$ is {\em $\phi$-flat\/} (or just $M$-flat, if
no ambiguity is likely to arise), if the functor
$$
M\text{-}\mathbf{Int.Mod}_\circ\to A\Mod \quad : \quad S\mapsto
S\otimes_MN
$$
is {\em exact}, in the sense that it sends exact sequences of pointed
integral $M$-modules, to exact sequences of $A$-modules. We say that
$N$ is {\em faithfully $\phi$-flat\/} if this functor is exact in the
above sense, and we have $S\otimes_MN=0$ if and only if $S=0$.
\end{definition}

\begin{remark}\label{rem_when-is-it-flat}
(i)\ \ 
Notice that the functor $S\mapsto S\otimes_MN$ of definition
\ref{def_M-flatness} is right exact in the categorical sense
({\em i.e.} it commutes with finite colimits), since it is a
right adjoint. However, even when $N$ is faithfully flat, this
functor is not always left exact in the categorical sense : it
does not commute with finite products, nor with equalizers, in
general.

(ii)\ \ 
Let $M$ and $S$ be as in definition \ref{def_M-flatness}, and
let $R$ be any non-zero commutative unital ring. Denote by
$\phi:M\to(R\La M\Ra,\cdot)$ the natural morphism of pointed
monoids. In light of \eqref{eq_yet-more-nonsense}, it is clear
that if $R\La S\Ra$ is a flat $R\La M\Ra$-module, then $S$ is
a flat pointed $M$-module, and the latter condition implies
that $R\La S\Ra$ is $\phi$-flat.

(iii)\ \
Let $P$ be a pointed monoid, $A$ a ring, $\phi:P\to(A,\cdot)$
a morphism of monoids, $f:A\to B$ a ring homomorphism.
If $A$ is $\phi$-flat and $f$ is flat, then $B$ is
$(f\circ\phi)$-flat. Conversely, if $B$ is $(f\circ\phi)$-flat
and $f$ is faithfully flat, then $A$ is $\phi$-flat : the
verifications are standard, and shall be left to the reader.
\end{remark}

\begin{lemma}\label{lem_faithful-phi-flat}
Let $M$ be a monoid, $A$ a ring, $\phi:M\to(A,\cdot)$ a morphism
of monoids, and assume that $\phi$ is local and $A$ is $\phi$-flat.
Then $A$ is faithfully $\phi$-flat.
\end{lemma}
\begin{proof} Let $S$ be an integral pointed $M$-module, and suppose
that $S\otimes_MA=\{0\}$; we have to show that $S=\{0\}$. Say that
$s\in S$, and let $Ms\subset S$ be the $M$-submodule generated by
$s$. Since $A$ is $\phi$-flat, it follows easily that
$Ms\otimes_MA=\{0\}$, hence we are reduced to the case where $S$ is
cyclic. By remark \ref{rem_integr-modules}(iii), we may then assume
that $S=M/I$ for some ideal $I\subset M$. It follows that
$S\otimes_MA=A/\phi(I)A$, so that $\phi(I)$ generates $A$. Since
$\phi$ is local, this implies that $I=M$, whence the contention.
\end{proof}

\begin{lemma}\label{lem_intersect-ideals}
Let $M$ be an integral pointed monoid, $I,J\subset M$ two ideals,
$A$ a ring, $\alpha:M\to(A,\cdot)$ a morphism of monoids, $N$ an
$\alpha$-flat $A$-module, and $S$ a flat $M$-module. Then :
$$
\begin{aligned}
IS\cap JS & =(I\cap J)S \\
\alpha(I)N\cap\alpha(J)N & =\alpha(I\cap J)N.
\end{aligned}
$$
\end{lemma}
\begin{proof} We consider the commutative ladder of pointed
$M$-modules, with exact rows and injective vertical arrows :
\set\begin{equation}\label{eq_ladder-of-mods}
\diagram 0 \ar[r] & I\cap J \ar[r] \ar[d] & I \ar[r] \ar[d] &
I/(I\cap J) \ar[r] \ar[d] & 0 \\
0 \ar[r] & J \ar[r] & M \ar[r] & M/(I\cup J) \ar[r] & 0
\enddiagram
\end{equation}
By assumption, the ladder of $A$-modules
$\eqref{eq_ladder-of-mods}\otimes_MN$ has still exact rows and
injective vertical arrows. Then, the snake lemma gives the following
short exact sequence involving the cokernels of the vertical
arrows :
$$
0 \to JN/(I\cap J)N\to N/IN\xrightarrow{p} N/(IN+JN)\to 0
$$
(where we have written $JN$ instead of $\alpha(J)N$, and likewise
for the other terms). However $\Ker\,p=JN/(IN\cap JN)$, whence the
second stated identity. The first stated identity can be deduced
from the second, by virtue of remark \ref{rem_when-is-it-flat}(ii).
\end{proof}

\begin{remark}\label{rem_whenever-phi}
By inspection of the proof, we see that the first identity
of lemma \ref{lem_intersect-ideals} holds, more generally,
whenever $\Z\La S\Ra$ is $\phi$-flat, where
$\phi:M\to\Z\La M\Ra$ is the natural morphism of pointed monoids.
\end{remark}

\begin{proposition}\label{prop_fla-criterion-point}
Let $M$ be a pointed integral monoid, $A$ a ring,
$\phi:M\to(A,\cdot)$ a morphism of monoids, $N$ an
$A$-module. Then we have :
\begin{enumerate}
\item
The following conditions are equivalent :
\begin{enumerate}
\item
$N$ is $\phi$-flat.
\item
$\Tor_i^{\Z\La M\Ra}(\Z\La T\Ra,N)=0$ for every $i>0$ and
every pointed integral $M$-module $T$.
\item
$\Tor_1^{\Z\La M\Ra}(\Z\La M/I\Ra,N)=0$ for every ideal
$I\subset M$.
\item
The natural map $I\otimes_MN\to N$ is injective for every
ideal $I\subset M$.
\end{enumerate}
\item
If moreover, $M^\sharp$ is finitely generated, then
the conditions {\em (a)-(d)} of {\em(i)} are equivalent
to either of the following two conditions :
\begin{enumerate}
\addenuii\addenuii\addenuii\addenuii
\item
$\Tor_1^{\Z\La M\Ra}(\Z\La M/\fp\Ra,N)=0$ for every prime
ideal $\fp\subset M$.
\item
The natural map $\fp\otimes_MN\to N$ is injective for every
prime ideal $\fp\subset M$.
\end{enumerate}
\end{enumerate}
\end{proposition}
\begin{proof} Clearly
(a)$\Rightarrow$(b)$\Rightarrow$(c)$\Rightarrow$(e).
Next, by considering the short exact sequence of pointed
integral $M$-modules $0\to I\to M\to M/I\to 0$ we easily
see that (c)$\Leftrightarrow$(d) and (e)$\Leftrightarrow$(f).

(c)$\Rightarrow$(a) : Let $\Sigma:=(0\to S'\to S\to S''\to 0)$
be a short exact sequence of pointed integral $M$-modules;
we need to show that the induced map $S'\otimes_MN\to S\otimes_MN$
is injective. Since the sequence $\Z\La\Sigma\Ra$ is still exact,
the long $\Tor$-exact sequence reduces to showing that
$\Tor_1^{\Z\La M\Ra}(\Z\La T\Ra,N)=0$ for every pointed integral
$M$-module $T$.
In view of remark \ref{rem_integr-modules}(i), we are easily
reduced to the case where $T$ is finitely generated; next,
using remark \ref{rem_integr-modules}(ii,iii), the long exact
$\Tor$-sequence, and an easy induction on the number of generators
of $T$, we may assume that $T=M/I$, whence the contention.

Lastly, if $M^\sharp$ is finitely generated, then remark
\ref{rem_integr-modules}(iv) shows that, in the foregoing
argument, we may further reduce to the case where $T=M/\fp$
for a prime ideal $\fp\subset M$; this shows that
(e)$\Rightarrow$(a).
\end{proof}

\begin{lemma}\label{lem_Lazard}
Let $M$ be a pointed monoid, $S$ a pointed $M$-module, and
suppose that the following conditions hold for $S$ :
\begin{itemize}
\item[(F1)]
If $s\in S$ and $a\in\Ann_M(s)$, then there exists $b\in\Ann_M(a)$
such that $s\in bS$.
\item[(F2)]
If $a_1,a_2\in M$ and $s_1,s_2\in S$ satisfy the identity
$a_1s_1=a_2s_2\neq 0$, then there exist $b_1,b_2\in M$ and
$t\in S$ such that $s_i=b_it$\ \ for $i=1,2$ and\ \ $a_1b_1=a_2b_2$.
\end{itemize}
Then the natural map $I\otimes_MS\to S$ is injective for every
ideal $I\subset M$.
\end{lemma}
\begin{proof} Let $I$ be an ideal, and suppose that two elements
$a_1\otimes s_1$ and $a_2\otimes s_2$ of $I\otimes_MS$ are
mapped to the same element of $S$.
If $a_is_i=0$ for $i=1,2$, then (F1) says that there exist
$b_1,b_2\in M$ and $t_1,t_2\in S$ such that $a_ib_i=0$ and
$s_i=b_it_i$ for $i=1,2$; thus
$a_i\otimes s_i=a_i\otimes b_it_i=a_ib_i\otimes s_i=0$ in
$I\otimes_MS$.
In case $a_is_i\neq 0$, pick $b_1,b_2\in M$ and $t\in S$
as in (F2); we conclude that
$a_1\otimes s_1=a_1\otimes b_1t=a_1b_1\otimes t=a_2b_2\otimes t=
a_2\otimes s_2$ in $I\otimes_MS$, whence the contention.
\end{proof}

\begin{theorem}\label{th_always-exact}
Let $M$ be an integral pointed monoid, $S$ a pointed $M$-module.
The following conditions are
equivalent :
\begin{enumerate}
\alphaenu
\item
$S$ is $M$-flat.
\item
For every morphism $M\to P$ of pointed monoids, $P\otimes_MS$
is $P$-flat.
\item
For every short exact sequence $\Sigma$ of integral pointed
$M$-modules, the sequence $\Sigma\otimes_MS$ is again short exact.
\item
Conditions {\em (F1)} and {\em (F2)} of lemma {\em\ref{lem_Lazard}}
hold for $S$.
\end{enumerate}
\end{theorem}
\begin{proof} Clearly (b)$\Rightarrow$(a)$\Rightarrow$(c).

(c)$\Rightarrow$(d): To show (F1), set $I:=Ma$, and denote
by $i:I\to M$ the inclusion; (c) implies that the induced
map $i\otimes_MS:I\otimes_MS\to S$ is injective. However,
we have a natural isomorphism $I\isom M/\Ann_M(a)$ of
$M$-modules (remark \ref{rem_integr-modules}(iii)), whence
an isomorphism : $I\otimes_MS\isom S/\Ann_M(a)S$, and
under this identification, $i\otimes_MS$ is induced
by the map $S\to S$ : $s\mapsto as$. Thus, multiplication
by $a$ maps the subset $S\setminus\!\Ann_M(a)S$ injectively
into itself, which is the claim.

For (F2), notice that $\Z\La S\Ra$ is $\phi$-flat under
condition (c), for $\phi:M\to\Z\La M\Ra$ the natural
morphism. Now, say that $a_1s_1=a_2s_2\neq 0$ in $S$;
set $I:=Ma_1$, $J:=Ma_2$; the assumption means that
$a_1s_1\in IS\cap JS$, in which case remark
\ref{rem_whenever-phi} shows that there exist $t\in S$
and $b_1,b_2\in M$ such that $a_1b_1=a_2b_2$, and
$a_1s_1=a_1b_1t$, hence $a_2s_2=a_2b_2t$. Since we have
seen that multiplication by $a_1$ maps
$S\setminus\!\Ann_M(a_1)S$ injectively into itself, we
deduce that $s_1=b_1t$, and likewise we get $s_2=b_2t$.

To prove that (d)$\Rightarrow$(b), we observe :

\begin{claim}\label{cl_loc-conn-cols}
Let $P$ be a pointed monoid, $\Lambda$ a small locally directed
category (see definition \ref{def_filtered-cols}(iv)), and
$S_\bullet:\Lambda\to P\Mod_\circ$ a functor, such that $S_\lambda$
fulfills conditions (F1) and (F2), for every $\lambda\in\Ob(\Lambda)$.
Then the colimit of $S_\bullet$ also fulfills conditions (F1) and (F2).
\end{claim}
\begin{pfclaim} In light of remark \ref{rem_cofinal}(ii), we may
assume that $\Lambda$ is either discrete or connected.
Suppose first that $\Lambda$ is connected; then remark
\ref{rem_point-and-complete}(ii) allows to check directly
that conditions (F1) and (F2) hold for the colimit of
$S_\bullet$, since they hold for every $S_\lambda$.
If $\Lambda$ is discrete, the assertion is that conditions
(F1) and (F2) are preserved by arbitrary (small) direct sums,
which we leave as an exercise for the reader.
\end{pfclaim}

To a given pointed $M$-module $S$, we attach the small category
$S^*$, such that :
$$
\Ob(S^*)=S\!\setminus\!\{0\}
\qquad \text{and} \qquad
\Hom_{S^*}(s',s)=\{a\in M~|~as=s'\}.
$$
The composition of morphisms is induced by the composition
law of $M$, in the obvious way. Notice that $S^*$ is locally
directed if and only if $S$ satisfies condition (F2).

We define a functor $F:S^*\to M\Mod_\circ$ as follows.
For every $s\in\Ob(S^*)$ we let $F(s):=M$, and for every
morphism $a:s'\to s$ we let $F(a):=a\cdot\one_M$.
We have a natural transformation $\tau:F\Rightarrow c_S$,
where $c_S:S^*\to M\Mod_\circ$ is the constant functor
associated with $S$; namely, for every $s\in\Ob(S^*)$,
we let $\tau_s:M\to S$ be the map given by the rule
$a\mapsto as$ for all $a\in M$. There follows a morphism
of pointed $M$-modules :
\set\begin{equation}\label{eq_loc_conn_colim}
\colim_{S^*}\,F\to S
\end{equation}

\begin{claim}\label{cl_F1-and-F2}
If $S$ fulfills conditions (F1), the map
\eqref{eq_loc_conn_colim} is an isomorphism.
\end{claim}
\begin{pfclaim} Indeed, we have a natural decomposition of
$S^*$ as coproduct of a family $(S^*_i~|~i\in I)$ of
connected subcategories (for some small set $I$ : see
remark \ref{rem_cofinal}(ii)); especially we have
$\Hom_{S^*}(s,s')=\emptyset$ if $s\in\Ob(S^*_i)$ and
$s'\in\Ob(S^*_j)$ for some $i\neq j$ in $I$.

For each $i\in I$, let $F_i:S^*_i\to M\Mod_\circ$ be the
restriction of $F$. There follows a natural isomorphism :
$$
\bigoplus_{i\in I}\colim_{S^*_i}F_i\isom\colim_{S^*}F.
$$
Since the colimit of $F_i$ commutes with the forgetful
functor to sets, an inspection of the definitions yields
the following explicit description of the colimit $T_i$
of $F_i$. Every element of $T_i$ is represented by some
pair $(s,a)$ where $s\in\Ob(S^*_i)\subset S\!\setminus\!\{0\}$
and $a\in M$; such pair is mapped to $as$ by
\eqref{eq_loc_conn_colim}, and two such pairs $(s,a)$, $(s',a')$
are identified in $T_i$ if there exists $b\in M$ such that 
$bs=s'$ and $ba'=a$.

Hence, denote by $S_i$ the image under \eqref{eq_loc_conn_colim}
of $T_i$; we deduce first, that $S_i\cap S_j=\{0\}$ if $i\neq j$.
Indeed, say that $t\in S_i\cap S_j$; by the foregoing, there exist
$s_i\in\Ob(S^*_i)$, $s_j\in\Ob(S^*_j)$, and $a_i,a_j\in M$, such
that $a_is_i=t=a_js_j$. If $t\neq 0$, we get morphisms $a_i:t\to s_i$
and $a_j:t\to s_j$ in $S^*$; say that $t\in S^*_k$ for some
$k\in I$; it then follows that $S^*_i=S^*_k=S^*_j$, a contradiction.
Next, it is clear that \eqref{eq_loc_conn_colim} is surjective.
It remains therefore only to show that each $T_i$ maps injectively
onto $S_i$. Hence, say that $(s_1,a_1)$ and $(s_2,a_2)$ represent
two elements of $T_i$ with $t:=a_1s_1=a_2s_2$. If $t\neq 0$,
we get, as before, morphisms $a_1:t\to s_1$ and $a_2:t\to s_2$
in $S^*_i$, and the two pairs are identified in $T_i$ to the pair
$(t,1)$.
Lastly, if $t=0$, condition (F1) yields $b\in M$ and $s'\in S$
such that $a_1b=0$ and $s_1=bs'$, whence a morphism $b:s_1\to s'$
in $S^*_i$, and $F_i(b)(a_1,s_1)=(0,s')$ which represents the zero
element of $T_i$. The same argument applies as well to $(a_2,s_2)$,
and the claim follows.
\end{pfclaim}

\begin{claim}\label{cl_inject-Fi}
Let $M\to P$ be any morphism of pointed monoids, and
$S$ a pointed $M$-module fulfilling conditions (F1) and (F2).
Then the natural map $I\otimes_MS\to P\otimes_MS$ is injective,
for every ideal $I\subset P$.
\end{claim}
\begin{pfclaim} From claim \ref{cl_F1-and-F2} we deduce that
$P\otimes_MS$ is the locally directed colimit of the functor
$P\otimes_MF$, and notice that the pointed $P$-module $P$ fulfills
conditions (F1) and (F2); by claim \ref{cl_loc-conn-cols} we
deduce that $P\otimes_MS$ also fulfills the same conditions,
so the claim follows from lemma \ref{lem_Lazard}.
\end{pfclaim}

After these preliminaries, suppose that conditions (F1) and (F2)
hold for $S$, and let $\Sigma:=(0\to T'\to T\to T''\to 0)$ be a
short exact sequence of pointed $M$-modules. We wish to show that
$\Sigma\otimes_MS$ is still short exact.
However, if $U''\subset T''$ is any $M$-submodule, let
$U\subset T$ be the preimage of $U''$, and notice that the
induced sequence $0\to T'\to U\to U''\to 0$ is still short exact.
Since a filtered colimit of short exact sequences is short exact,
remark \ref{rem_integr-modules}(i) allows to reduce to the case
where $T''$ is finitely generated.

We shall argue by induction on the number $n$ of generators
of $T''$. Hence, suppose first that $T''$ is cyclic, and
let $t\in T$ be any element whose image in $T''$ is a
generator. Set $C:=Mt\subset T$, and let $C'\subset T'$
be the preimage of $C$. We obtain a cocartesian (and cartesian)
diagram of pointed $M$-modules :
$$
\cD \quad:\quad
{\diagram C' \ar[r] \ar[d] & C \ar[d] \\
          T' \ar[r] & T.
\enddiagram}\qquad\qquad
$$
The induced diagram $\cD\otimes_MS$ is still cocartesian, hence
the same holds for the diagram of sets underlying $\cD\otimes_MS$
(remark \ref{rem_point-and-complete}(ii)). Especially, if
the induced map $C'\otimes_MS\to C\otimes_MS$ is injective,
the same will hold for the map $T'\otimes_MS\to T\otimes_MS$.
We may thus replace $T'$ and $T$ by respectively $C'$ and $C$,
which allows to assume that also $T$ is cyclic. In this case,
pick a generator $u\in T$; we claim that there exists a unique
multiplication law $\mu_T$ on $T$, such that the surjection
$p:M\to T$ : $a\mapsto au$ is a morphism of pointed monoids.
Indeed, for every $t,t'\in T$, write $t=au$ for some $a\in M$,
and set $\mu_T(t,t'):=at'$. Using the linearity of $p$ we
easily check that $\mu_T(t,t')$ does not depend on the choice
of $a$, and the resulting composition law $\mu_T$ is commutative
and associative. Then $T'$ is an ideal of $T$, so claim
\ref{cl_inject-Fi} tells us that the map
$T'\otimes_MS\to T\otimes_MS$ is injective, as required.

Lastly, suppose that $n>1$, and the assertion is already known
whenever $T''$ is generated by at most $n-1$ elements.
Let $U''\subset T''$ be a pointed $M$-submodule, such that
$U''$ is generated by at most $n-1$ elements, and $T''/U''$
is cyclic. Denote by $U\subset T$ the preimage of $U''$;
we deduce short exact sequences
$\Sigma':=(0\to T'\to U\to U''\to 0)$ and
$\Sigma'':=(0\to U\to T\to T''/U''\to 0)$, and by inductive
assumption, both $\Sigma'\otimes_MS$ and $\Sigma''\otimes_MS$
are short exact. Therefore the natural map
$T'\otimes_MS\to T\otimes_MS$ is the composition of two
injective maps, hence it is injective, as stated.
\end{proof}

\begin{remark}\label{rem_when-is-it-flat-II}
In the situation of remark \ref{rem_when-is-it-flat}(ii),
suppose that $M$ is pointed integral. Then theorem
\ref{th_always-exact} implies that $S$ is a flat pointed
$M$-module if and only if $R\La S\Ra$ is $\phi$-flat.
\end{remark}

\begin{corollary}\label{cor_first-corollo}
Let $M$ be an integral pointed monoid, $S$ a pointed $M$-module.
Then
\begin{enumerate}
\item
The following conditions are equivalent :
\begin{enumerate}
\item
$S$ is $M$-flat.
\item
For every ideal $I\subset M$, the induced map $I\otimes_MS\to S$
is injective.
\end{enumerate}
\item
If moreover $M^\sharp$ is finitely generated, then these
conditions are equivalent to :
\begin{enumerate}
\addenuii
\addenuii
\item
For every prime ideal $\fp\subset M$, the induced map
$\fp\otimes_MS\to S$ is injective.
\end{enumerate}
\end{enumerate}
\end{corollary}
\begin{proof}(i): Indeed, by remark \ref{rem_when-is-it-flat-II}
and proposition \ref{prop_fla-criterion-point}(i) (together with
\eqref{eq_yet-more-nonsense}), both (a) and (b) hold if and only
if $\Z\La S\Ra$ is $\phi$-flat, for $\phi:M\to\Z\La M\Ra$ the
natural morphism.

(ii): This follows likewise from proposition
\ref{prop_fla-criterion-point}(ii).
\end{proof}

\begin{corollary}\label{cor_yet-another-flat}
Let $\phi:M\to N$ be a morphism of pointed monoids,
$G\subset M^\times$, $H\subset N^\times$ two subgroups
such that $\phi(G)\subset H$, and denote by
$\bar\phi:M/G\to N/H$ the induced morphism.
Let also $S$ be any $N$-module. We have :
\begin{enumerate}
\item
If $S_{(\phi)}$ is a flat $M$-module, then $S/H_{(\bar\phi)}$
is a flat $M/G$-module (notation of \eqref{subsec_tens-restr}).
\item
If moreover, $M$ is a pointed integral monoid and $S$
is a pointed integral $H$-module, then also the converse
of {\em(i)} holds.
\end{enumerate}
\end{corollary}
\begin{proof}(i): We have a natural isomorphism
$$
(S/H)_{(\bar\phi)}\isom
N/H\otimes_{N/\phi G}S_{(\phi)}/\phi G.
$$
However, by example \ref{ex_flat-quots}(ii), $N/H$
is a flat $N/\phi G$-module, and $S_{(\phi)}/\phi G$
is a flat $M/G$-module, whence the contention.

(ii): By theorem \ref{th_always-exact}, it suffices to check
that conditions (F1) and (F2) of lemma \ref{lem_Lazard} hold
in $S_{(\phi)}$, and notice that, by the same token, both
conditions hold for the $M/G$-module $S/H_{(\bar\phi)}=S_{(\phi)}/G$,
since $M/G$ is pointed integral (example \ref{ex_flat-quots}(v)).

Hence, say that $\phi(a)s=0$ for some $a\in M$ and $s\in S$;
we may then find $b\in M$, $t'\in S$ and $g\in G$ such that
$\phi(gb)t=s$ and $ab=0$. Setting $t:=\phi(g)t'$, we deduce
that (F1) holds.

Next, say that $\phi(a_1)s_1=\phi(a_2)s_2\neq 0$ for some
$a_1,a_2\in M$ and $s_1,s_2\in S$; then we may find
$g\in G$, $h_1,h_2\in H$, $b_1,b_2\in M$ and $t'\in S$ such that
$ga_1b_1=a_2b_2$ and
\set\begin{equation}\label{eq_wasp}
\phi(b_i)h_it'=s_i
\qquad
\text{for $i=1,2$}.
\end{equation}
After replacing $b_1$ by $gb_1$ and $h_1$ by $h_1\phi(g^{-1})$,
we reduce to the case where $a_1b_1=a_2b_2$. From \eqref{eq_wasp}
we deduce that
$$
\phi(a_1b_1)h_1t'=\phi(a_1)s_1=\phi(a_2)s_2=\phi(a_2b_2)h_2t'=
\phi(a_1b_1)h_2t'
$$
whence $h_1=h_2$, since $S$ is a pointed integral $H$-module.
Setting $t:=h_1t'$, we see that (F2) holds.
\end{proof}

\begin{corollary}\label{cor_flat-face}
Let $\phi:M\to N$ be a flat morphism of pointed monoids, with
$M$ pointed integral, and let $\fp\subset N$ be any prime ideal.
Then the morphism $M/\phi^{-1}\fp\to N/\fp$ induced by $\phi$ is
also flat.
\end{corollary}
\begin{proof} Let $F:=N\!\setminus\!\fp$; by theorem
\ref{th_always-exact}, it suffices to check that conditions
(F1) and (F2) of lemma \ref{lem_Lazard} hold for the
$\phi^{-1}F$-module $F$. However, since $0\notin F$,
condition (F1) holds trivially. Moreover, by assumption these
two conditions hold for the $M$-module $N$; hence, say that
$\phi(a_1)\cdot s_1=\phi(a_2)\cdot s_2$ in $F$, for some
$a_1,a_2\in\phi^{-1}F$ and $s_1,s_2\in F$. It follows that
there exist $b_1,b_2\in M$ and $t\in N$ such that $a_1b_1=a_2b_2$
in $M$, and $\phi(b_i)\cdot t=s_i$ for $i=1,2$. Since $F$ is
a face, this implies that $\phi(b_1),\phi(b_2),t\in F$, so
(F2) holds for $F$, as stated.
\end{proof}

Another corollary is the following analogue of a well
known criterion due to Lazard.

\begin{proposition}\label{prop_Lazard}
Let $M$ be an integral pointed monoid, $S$ an integral pointed
$M$-module, $R$ a non-zero commutative ring with unit. The
following conditions are equivalent :
\begin{enumerate}
\alphaenu
\item
$S$ is $M$-flat.
\item
$S$ is the colimit of a filtered system of free pointed
$M$-modules (see remark {\em\ref{rem_point-and-complete}(ii)}).
\item
$R\La S\Ra$ is a flat $R\La M\Ra$-module.
\end{enumerate}
\end{proposition}
\begin{proof} Obviously (b)$\Rightarrow$(a) and
(b)$\Rightarrow$(c).

(c)$\Rightarrow$(a) has already been observed in remark
\ref{rem_when-is-it-flat}(ii).

(a)$\Rightarrow$(b): It suffices to prove that if $S$
is flat, the category $S^*$ attached to $S$ as in the proof
of theorem \ref{th_always-exact}, is pseudo-filtered.
However, a simple inspection of the construction shows
that $S^*$ is pseudo-filtered if and only if $S$ satisfies
both condition (F2) of lemma \ref{lem_Lazard}, and
the following further condition. For every $a,b\in M$ and
$s\in S$ such that $as=bs\neq 0$, there exist $t\in S$ and
$c\in M$ such that $ac=bc$ and $ct=s$. This condition is
satisfied by every integral pointed $M$-module, whence the
contention.
\end{proof}

\sset\subsubsection{}\label{subsec_cart-mons}
Consider now, a cartesian diagram of integral pointed monoids :
$$
\cD(P_0,I,P_1)
\quad : \quad
{\diagram P_0 \ar[r] \ar[d] & P_1 \ar[d] \\
          P_2 \ar[r] & P_3
\enddiagram}
$$
where $P_2$ (resp. $P_3$) is a quotient $P_0/I$ (resp. $P_1/IP_1$)
for some ideal $I\subset P_0$, and the vertical arrows of
$\cD(P_0,I,P_1)$ are the natural surjections.
In this situation, it is easily seen that the induced map
$I\to IP_1$ is bijective; especially, $I$ is both a $P_0$-module
and a $P_1$-module.
Let $\phi_i:P_i\to\Z\La P_i\Ra$ be the units of adjunction,
for $i=0,1,2$. Let also $M$ be any $\Z\La P_0\Ra$-module.

\begin{lemma}\label{lem_not-flat-desc}
In the situation of \eqref{subsec_cart-mons}, we have :
\begin{enumerate}
\item
Let $J\subset P_0$ be any ideal. Then $S:=P_0/J$ admits a
three-step filtration
$$
0\subset\Fil_0S\subset\Fil_1S\subset\Fil_2S=S
$$
such that $\Fil_0S$ and $\gr_2S$ are $P_2$-modules, and
$\gr_1S$ is a $P_1$-module.
\item
The following conditions are equivalent :
\begin{enumerate}
\item
$M$ is $\phi_0$-flat.
\item
$M\otimes_{P_0}P_i$ is $\phi_i$-flat and
$\Tor_1^{\Z\La P_0\Ra}(M,\Z\La P_i\Ra)=0$, for $i=1,2$.
\end{enumerate}
\item
The following conditions are equivalent :
\begin{enumerate}
\item
$\Tor_1^{\Z\La P_0\Ra}(M,\Z\La P_i\Ra)=0$ for $i=1,2,3$.
\item
$\Tor_1^{\Z\La P_i\Ra}(M\otimes_{P_0}P_i,\Z\La P_3\Ra)=0$
for $i=1,2$.
\end{enumerate}
\item
Suppose moreover, that $P_3\neq 0$ is a free pointed
$P_2$-module. Then the $\Z\La P_0\Ra$-module $M$ is
$\phi_0$-flat if and only if the $\Z\La P_1\Ra$-module
$M\otimes_{P_0}P_1$ is $\phi_1$-flat.
\end{enumerate}
\end{lemma}
\begin{proof}(i): Define
$\Fil_1S:=\Ker(S\to P_0/(I\cup J))=(I\cup J)/J$.
Then it is already clear that $S/\Fil_1S$ is a $P_2$-module.
Next, let $\Fil_0S:=\Ker(\Fil_1S\to(I\cup JP_1)/JP_1)=JP_1/J$.
Since $IP_1=I$, we see that $\Fil_0S$ is a $P_2$-module,
and $(I\cup JP_1)/JP_1$ is a $P_1$-module.

(ii): Clearly (a)$\Rightarrow$(b), hence suppose that (b) holds,
and let us prove (a). By proposition \ref{prop_fla-criterion-point}(i),
it suffices to show that $\Tor_1^{\Z\La P_0\Ra}(M,\Z\La P_0/I\Ra)=0$
for every ideal $J\subset P_0$. In view of (i), we are then
reduced to showing that $\Tor_1^{\Z\La P_0\Ra}(M,\Z\La S\Ra)=0$,
whenever $S$ is a $P_i$-module, for $i=1,2$. However, for such
$P_i$-module $S$, we have a base change spectral sequence
$$
E^2_{pq}:
\Tor_p^{\Z\La P_i\Ra}(\Tor_q^{\Z\La P_0\Ra}(M,\Z\La P_i\Ra),\Z\La S\Ra)
\Rightarrow\Tor_{p+q}^{\Z\La P_0\Ra}(M,\Z\La S\Ra).
$$
Under assumption (b), we deduce :
$\Tor_1^{\Z\La P_0\Ra}(M,\Z\La S\Ra)=
\Tor_1^{\Z\La P_i\Ra}(M\otimes_{P_0}P_i,\Z\La S\Ra)=0$.

(iii): Notice that the induced diagram of rings
$\Z\La\cD(P_0,I,P_1)\Ra$ is still cartesian. Then, this is a
special case of \cite[Lemma 3.4.15]{Ga-Ra}.

(iv): Suppose that $M\otimes_{P_0}P_1$ is $\phi_1$-flat, and $P_3$
is a free pointed $P_2$-module, say $P_3\simeq P_2^{(\Lambda)_\circ}$,
for some set $\Lambda\neq\emptyset$; then the $\Z\La P_2\Ra$-module
$\Z\La P_3\Ra$ is isomorphic to $\Z\La P_2\Ra^{(\Lambda)}$,
especially it is faithfully flat, and we deduce that
$\Tor^{\Z\La P_0\Ra}_1(M,\Z\La P_i\Ra)=0$ for $i=1,2$, by (iii).
On the other hand, $M\otimes_{P_0}P_3$ is $\phi_3$-flat,
so it also follows that $M\otimes_{P_0}P_2$ is $\phi_2$-flat,
by proposition \ref{prop_Lazard}. Summing up, this shows that
$M$ fulfills condition (ii.b), hence also (ii.a), as sought.
\end{proof}

\sset\subsubsection{}\label{subsec_not-flat-desc}
Let now $P\to Q$ be an injective morphism of integral pointed
monoids, and suppose that $P^\sharp$ and $Q^\sharp$ are finitely
generated monoids, and $Q$ is a finitely generated $P$-module.
Denote $\phi_P:P\to\Z\La P\Ra$ and $\phi_Q:Q\to\Z\La Q\Ra$ the
usual units of adjunction.

\begin{theorem} In the situation of \eqref{subsec_not-flat-desc},
let $M$ be a $\Z\La P\Ra$-module, and suppose that the
$\Z\La Q\Ra$-module $M\otimes_PQ$ is $\phi_Q$-flat. Then
$M$ is $\phi_P$-flat.
\end{theorem}
\begin{proof} Using lemma \ref{lem_not-flat-desc}(iv), and an
easy induction, it suffices to show that there exists a finite
chain
\set\begin{equation}\label{eq_chain-pinteg}
P=Q_0\subset Q_1\subset\cdots\subset Q_n=Q
\end{equation}
of inclusions of integral pointed monoids, and for every
$j=0,\dots,n-1$ an ideal $I_j\subset Q_j$, and a cartesian
diagram of integral pointed monoids $\cD(Q_j,I_j,Q_{j+1})$
as in \eqref{subsec_cart-mons}, and such that $Q_{j+1}/I_j\neq 0$
is a free pointed $Q_j/I_j$-module. (Notice that each $Q_i^\sharp$
is a quotient of $Q_i/P^\times$, and the latter is a submodule
of $Q/P^\times$, hence $Q_i^\sharp$ is still a finitely
generated monoid, by proposition \ref{prop_ideals-in-fg-mon}(i).)
If $P=Q$, there is nothing to prove; so we may assume that $P$
is strictly contained in $Q$, and -- invoking again proposition
\ref{prop_ideals-in-fg-mon}(i) -- we further reduce to showing
that there exist a monoid $Q_1\subset Q$ strictly containing
$P$, and an ideal $I\subset P$, such that the diagram
$\cD(P,I,Q_1)$ fulfills the above conditions.

Suppose first that the support of $Q/P$ contains only the
maximal ideal $\fm_P$ of $P$, and let $x_1,\dots,x_r$ be
a finite system of generators for $\fm_P$ (proposition
\ref{prop_ideals-in-fg-mon}(ii)). For each $i=1,\dots,r$,
the localization $(Q/P)_{x_i}$ is a $P_{x_i}$-module with
empty support, hence $(Q/P)_{x_i}=0$ (remark
\ref{rem_support-mods-mon}(v)). It follows that every
element of $Q/P$ is annihilated by some power of $x_i$,
and since $Q/P$ is finitely generated, we may find
$N\in\N$ large enough, such that $x_i^NQ/P=0$ for
$i=1,\dots,r$. After replacing $N$ by some possibly larger
integer, we get $\fm_P^N\cdot Q/P=0$, and we may assume
that $N$ is the least integer with this property. If $N=0$,
there is nothing to prove; hence suppose that $N>0$, and
set $Q_1:=P\cup\fm_P^{N-1}Q$. Notice that $Q_1$ is a monoid,
and $Q_1/P\neq 0$ is annihilated by $\fm_P$. Especially,
$\fm_PQ_1=\fm_P$. Moreover, the induced map
$P/\fm_P\to Q_1/\fm_P$ is injective and $Q_1$ is an integral
pointed module, therefore the group $P^\times$ acts freely
on $Q_1\!\setminus\!\fm_P$, {\em i.e.} $Q_1/\fm_P$ is a free
pointed $P/\fm_P$-module.
It follows easily that if we take $I:=\fm_P$, we do obtain
a diagram $\cD(P,\fm_P,Q_1)$ with the sought properties, in
this case.

For the general case, let $\fp\subset P$ be a minimal
element of $\Supp\,Q/P$ (for the ordering given by inclusion).
Then the induced morphism $P_\fp\to Q_\fp$ still satisfies
the conditions of \eqref{subsec_not-flat-desc} (lemma
\ref{lem_face}(iv)). Moreover, $\Supp\,Q_\fp/P_\fp=\{\fp P_\fp\}$
by remark \ref{rem_support-mods-mon}(vii). By the previous case,
we deduce that there exists a chain of inclusions of integral
pointed monoids $P_\fp\subset Q'_1\subset Q_\fp$, such that
the resulting diagram $\cD(P_\fp,\fp P_\fp,Q'_1)$ is cartesian,
and $Q'_1/\fp P_\fp\neq 0$ is a free pointed $P_\fp/\fp P_\fp$-module.
Let $\bar e_1,\dots,\bar e_d$ be a basis of the latter
$P_\fp/\fp P_\fp$-module, with $\bar e_1=1$. Hence,
$\bar e_i\in Q_\fp\setminus\fp P_\fp$ for every $i=1,\dots,d$,
and after multiplying $\bar e_2,\dots,\bar e_d$ by a suitable
element of $P\setminus\fp$, we may assume that each
$\bar e_i$ is the image in $Q_\fp$ of an element $e_i\in Q$.
Moreover, for every $i,j\leq d$, either $\bar e_i\bar e_j=0$,
or else there exist $a_{ij}\in P_\fp$ and $k(i,j)\leq d$ such
that $\bar e_i\bar e_j=a_{ij}\bar e_{k(i,j)}$.
Furthermore, fix a system of generators $x_1,\dots,x_r$ for
$\fp$; then, for every $i\leq r$ and every $j\leq d$ we have
$x_i\bar e_j\in\fp P_\fp$. Again, after multiplying
$\bar e_2,\dots,\bar e_d$ by some $c\in P\!\setminus\!\fp$,
we may assume that $a_{ij}\in P$ for every $i,j\leq d$, and
moreover that $x_i\bar e_j$ lies in the image of $\fp$ for
every $i\leq r$ and $j\leq d$.

And if we multiply yet again $e_2,\dots,e_d$ by a suitable
element of $P\!\setminus\!\fp$, we may finally reach a
system of elements $e_1,\dots,e_d\in Q$ such that $e_1=1$ and :
\begin{itemize}
\item
For every $i,j\leq d$, we have either $e_ie_j=0$ or else
$e_ie_j=a_{ij}e_{k(i,j)}$.
\item
$x_ie_j\in\fp$ for every $i\leq r$ and $j\leq d$.
\end{itemize}
Clearly these elements span a $P$-module $Q_1$ which is
a monoid containing $P$ and contained in $Q$; moreover,
by construction we have $\fp Q_1=\fp$, hence the resulting
diagram $\cD(P,\fp,Q_1)$ is cartesian. Notice also that
$(Q_1/\fp)_\fp\simeq Q'_1/\fp P_\fp$, and that $P/\fp=P'_{\!\circ}$,
where $P':=P\!\setminus\!\fp$ is an integral (non-pointed)
monoid. To conclude it suffices now to apply the following

\begin{claim} Let $P'$ be an integral non-pointed monoid, $S$
a pointed $P'_{\!\circ}$-module, and $\underline e:=(e_1,\dots,e_d)$
a system of generators for $S$. Suppose that
$S\otimes_{P'_{\!\circ}}P^{\prime\gp}_{\!\circ}$ is a free pointed
$P^{\prime\gp}_{\!\circ}$-module, and the image of $\underline e$
is a basis for this module. Then $S$ is a free pointed
$P'_{\!\circ}$-module, with basis $\underline e$.
\end{claim}
\begin{pfclaim}[] If $\underline e$ is not a basis, we have
a relation in $S$ of the type $a_1e_1=a_2e_2$, for some
$a_1,a_2\in P'$. This relation must persist in
$S\otimes_{P'_{\!\circ}}P^{\prime\gp}_{\!\circ}$, and implies
that $a_1=a_2=0$ in $P^{\prime\gp}_{\!\circ}$. However, under
the stated assumptions the localization map
$P'_{\!\circ}\to P_{\!\circ}^{\prime\gp}$ is injective,
a contradiction.
\end{pfclaim}
\end{proof}

\subsection{Integral monoids}\label{subsec_intg-monnd}
We begin presently the study of a special class of monoids,
the integral non-pointed monoids, and the subclass of
saturated monoids (see definition \ref{def_exact-phi}(iii)).
Later, we shall complement this section with further results
on fine monoids (see sections \ref{sec_fine-satu} and
\ref{sec_log-regular}). {\em Throughout this section,
all the monoids under consideration shall be non-pointed.}

\begin{definition} Let $\phi:M\to N$ be a morphism of monoids.
\begin{enumerate}
\item
$\phi$ is said to be {\em integral\/} if, for any integral monoid
$M'$, and any morphism $M\to M'$, the push-out $N\otimes_MM'$ is
integral.
\item
$\phi$ is said to be {\em strongly flat\/} (resp.
{\em strongly faithfully flat}) if the induced morphism
$\Z[\phi]:\Z[M]\to\Z[N]$ is flat (resp. faithfully flat).
\end{enumerate}
\end{definition}

\begin{lemma}\label{lem_obvious-int}
Let $f:M\to N$ and $g:N\to P$ be two morphisms of monoids.
\begin{enumerate}
\item
If $f$ and $g$ are integral (resp. strongly flat), the same
holds for $g\circ f$.
\item
If $f$ is integral (resp. strongly flat), and $M\to M'$ is any other
morphism, then the morphism $\one_{M'}\otimes_Mf:M'\to M'\otimes_MN$
is integral (resp. strongly flat).
\item
If $f$ is integral, and $S\subset M$ and $T\subset N$ are
any two submonoids such that $f(S)\subset T$, then the induced
morphism $S^{-1}M\to T^{-1}N$ is integral.
\item
If $S\subset M$ is any submonoid, the natural map $M\to S^{-1}M$ is
strongly flat.
\item
If $f$ is integral, then the same holds for $f^\intg$, and the
natural map $M^\intg\otimes_MN\to N^\intg$ is an isomorphism.
\item
The unit of adjunction $M\to M^\intg$ is an integral morphism.
\end{enumerate}
\end{lemma}
\begin{proof} (i) is obvious. Assertion (ii) for integral
maps is likewise clear, and (ii) for strongly flat morphisms
follows from \eqref{eq_push-out-tensor}. Assertion (iv)
follows immediately from \eqref{eq_loc-monds}.

(iii): Let $S^{-1}M\to M'$ be any map of integral
monoids; in view of lemma \ref{lem_localize}, we have
$P:=M'\otimes_{S^{-1}M}T^{-1}N\simeq T^{-1}(M'\otimes_MN)$,
hence $P$ is integral, which is the claim.

(v): The second assertion follows easily by comparing the
universal properties (details left to the reader); using
this and and (ii), we deduce that $f^\intg$ is integral.

(vi) is left to the reader.
\end{proof}

\begin{theorem}\label{th_flat-crit-for-mnds}
Let $\phi:M\to N$ be a morphism of integral monoids.
Consider the following conditions :
\begin{enumerate}
\alphaenu
\item
$\phi$ is integral.
\item
$\phi$ is flat (see remark {\em\ref{rem_flatness}(vi)}).
\item
$\phi$ is flat and injective.
\item
$\phi$ is strongly flat.
\item
For every field $k$, the induced map $k[\phi]:k[M]\to k[N]$ is flat.
\end{enumerate}
Then : {\em (e) $\Leftrightarrow$ (d) $\Leftrightarrow$ (c)
$\Rightarrow$ (b)
 $\Leftrightarrow$ (a)}.
\end{theorem}
\begin{proof} This result appears in \cite[Prop.4.1(1)]{Ka},
with a different proof.

(a)$\Rightarrow$(b): Let $I\subset M$ be any ideal; we consider
the $M$-module :
$$
\sR(M,I):=\bigoplus_{n\in\N}I^n
$$
where $I^n$ denotes, for each $n\in\N$, the $n$-th power of
$I$ in the monoid $(\cP(M),\cdot)$ of \eqref{subsec_toric}.
Then there exists an obvious multiplication law on $\sR(M,I)$,
such that the latter is a $\N$-graded integral monoid, and
the inclusion $M\to\sR(M,I)$ in degree zero is a morphism
of monoids. We call $\sR(M,I)$ the {\em Rees monoid\/}
associated with $M$ and $I$.

Denote also by $j:\sR(M,I)\to M\times\N$ the
natural inclusion map. By assumption, the monoid
$\sR(M,I)\otimes_MN$ is integral. However, the natural map
$\sR(M,I)\otimes_MN\to(\sR(M,I)\otimes_MN)^\gp$ factors through
$j\otimes_MN$. The latter means that, for every $n\in\N$,
the induced map $I^n\otimes_MN\to N$ is injective. Then
the assertion follows from proposition
\ref{prop_fla-criterion-point} and remark
\ref{rem_point-unpointed-topos}(ii).

(b)$\Rightarrow$(a): Let $M\to M'$ be any morphism of
integral monoids; we need to show that the natural map
$M'\otimes_MN\to M^{\prime\gp}\otimes_{M^\gp}N^\gp$ is
injective (see remark \ref{rem_tens-is-pushout}(i)).
The latter factors through the morphism
$M'\otimes_MN\to M^{\prime\gp}\otimes_MN$, which is
injective by theorem \ref{th_always-exact} and remark
\ref{rem_point-unpointed-topos}(ii).  We are thus reduced
to proving the injectivity of the natural map :
$$
M^{\prime\gp}\otimes_{M^\gp}(M^\gp\otimes_MN)\to
M^{\prime\gp}\otimes_{M^\gp}N^\gp.
$$
By comparing the respective universal properties, it is easily
seen that $M^\gp\otimes_MN$ is the localization $\phi(M)^{-1}N$,
which of course injects into $N^\gp$. Then the contention
follows from the following general :

\begin{claim} Let $G$ be a group, $T\to T'$ an injective
morphism of monoids, $G\to P$ a morphism of monoids.
Then the natural map $P\otimes_GT\to P\otimes_GT'$ is
injective.
\end{claim}
\begin{pfclaim} This follows easily from remark
\ref{rem_tens-is-pushout}(i) and lemma
\ref{lem_special-p-out}(ii).
\end{pfclaim}

(d)$\Rightarrow$(e) and (c)$\Rightarrow$(b) are trivial.

(e)$\Rightarrow$(c): The flatness of $\phi$ has already
been noticed in remark \ref{rem_when-is-it-flat}(ii).
To show that $\phi$ is injective, let $a_1,a_2\in M$ and
let $k$ be any field.
Under assumption (e), the image in $k[N]$ of the annihilator
$\Ann_{k[M]}(a_1-a_2)$ generates the ideal
$\Ann_{k[N]}(\phi(a_1)-\phi(a_2))$. Set $b:=a_1a_2^{-1}$;
it follows that $\Ann_{k[M^\gp]}(1-b)$ generates
$\Ann_{k[N^\gp]}(1-\phi(b))$. However, one checks easily
that the annihilator of $1-b$ in $k[M^\gp]$ is either $0$
if $b$ is not a torsion element of $M^\gp$, or else is
generated (as an ideal) by $1+b+\cdots+b^{n-1}$, where $n$
is the order of $b$ in the group $M^\gp$. Now, if
$\phi(a_1)=\phi(a_2)$, we have $\phi(b)=1$, hence the
annihilator of $1-b$ cannot be $0$, and in fact
$k[N^\gp]=\Ann_{k[N^\gp]}(1-\phi(b))=nk[N^\gp]$. Since $k$
is arbitrary, it follows that $n=1$, {\em i.e.} $a_1=a_2$.

(c)$\Rightarrow$(d): since $\phi$ is injective, $N_\circ$ is
an integral pointed $M_\circ$-module, so the assertion is
a special case of proposition \ref{prop_Lazard}.
\end{proof}

\begin{remark}\label{rem_why-not}
(i)\ \
Let $G$ be a group; then every morphism of monoids $M\to G$
is integral. Indeed, lemma \ref{lem_obvious-int}(i,vi) reduces
the assertion to the case where $M$ is integral, in which
case it is an immediate consequence of theorem
\ref{th_flat-crit-for-mnds} and example \ref{ex_flat-quots}(vi).

(ii)\ \
Let $M$ be an integral monoid, and $S$ a flat pointed
$M_\circ$-module. From theorem \ref{th_always-exact} we see
that $T:=S\!\setminus\!\{0\}$ is an $M$-submodule of $S$,
hence $S=T_\circ$.

(iii)\ \
Let $\phi:M\to N$ be a morphism of integral monoids, and set
$\Gamma:=\Coker\,\phi^\gp$; then the natural map $\pi:N\to\Gamma$
defines a grading on $N$ (see definition \ref{def_grad-monoids}(i)),
which we call the {\em $\phi$-grading}. As usual, we shall
write $N_\gamma$ instead of $\pi^{-1}(\gamma)$, for every
$\gamma\in\Gamma$. We shall use the additive notation for
the composition law of $\Gamma$; especially, the neutral
element shall be denoted by $0$. Clearly $\phi$ factors
through a morphism of monoids $M\to N_0$, and each graded
summand $N_\gamma$ is naturally a $M$-module.

(iv)\ \
With the notation of (iii), we claim that the induced morphism
$\bar\phi:\phi(M)\to N$ is flat (hence strongly flat, by theorem
\ref{th_flat-crit-for-mnds}), if and only if $N_\gamma$ is a
filtered union of cyclic $M$-modules, for every $\gamma\in\Gamma$.
Indeed, notice that a cyclic $M$-submodule of $N_\gamma$
is a free $\phi(M)$-module of rank one (since $N$ is integral),
hence the condition implies that $N_\gamma$ is a flat
$\phi(M)$-module, hence $\bar\phi$ flat. Conversely, suppose that
$\bar\phi$ is flat, and let $n_1,n_2\in N_\gamma$ (for some
$\gamma\in\Gamma$); this means that there exist $a_1,a_2\in M$
such that $\phi(a_1)n_1=\phi(a_2)n_2$ in $N$. Then, condition
(F2) of theorem \ref{th_always-exact} says that there exist
$n'\in N$ and $b_1,b_2\in M$ such that $n_i=\phi(b_i)n'$ for
$i=1,2$; especially, $n'\in N_\gamma$, which shows that $N_\gamma$
is a filtered union of cyclic $M$-modules.

(v)\ \
With the notation of (iii), notice as well, that a morphism
$\phi:M\to N$ of integral monoids is exact (see definition
\ref{def_exact-phi}(i)) if and only if $\Ker\,\phi^\gp\subset M$
and $\phi$ induces an isomorphism $M/\Ker\,\phi^\gp\isom N_0$.
(Details left to the reader.)
\end{remark}

\begin{theorem}\label{th_Gruson}
Let $M\to N$ be a finite, injective morphism of integral monoids,
and $S$ a pointed $M_\circ$-module. Then $S$ is $M_\circ$-flat if
and only if $N_\circ\otimes_{M_\circ}S$ is $N_\circ$-flat.
\end{theorem}
\begin{proof} In light of theorem \ref{th_always-exact}, we may
assume that $N_\circ\otimes_{M_\circ}S$ is $N_\circ$-flat, and
we shall show that $S$ is $M_\circ$-flat. To this aim, let
$1$ denote the trivial monoid (the initial and final object
in the category of monoids); any pointed $M_\circ$-module $X$
is a pointed $1_\circ$-module by restriction of scalars, and
if $X$ and $Y$ are any two pointed $M_\circ$-modules, we define
a $M_\circ$-module structure on $X\otimes_{1_\circ}Y$ by the
rule
$$
a\cdot(x\otimes y):=ax\otimes ay
\qquad
\text{for every $a\in M_\circ$, $x\in X$ and $y\in Y$}.
$$
With this notation, we remark :

\begin{claim}\label{cl_prod_with-G}
Let $\phi:M\to G$ be a morphism of monoids, where $G$ is a
group and $M$ is integral. Then a pointed $M_\circ$-module
$S$ is $M_\circ$-flat if and only if the same holds for the
$M_\circ$-module $S\otimes_{1_\circ}G_\circ$.
\end{claim}
\begin{pfclaim} Suppose that $S$ is $M_\circ$-flat; then
$S=T_\circ$ for some $M$-module $T$ (remark \ref{rem_why-not}(ii)),
and it is easily seen that
$S\otimes_{1_\circ}G_\circ=(T\times G)_\circ$. By
theorem \ref{th_always-exact}, it suffices to check that
conditions (F1) and (F2) of lemma \ref{lem_Lazard} hold
for $(T\times G)_\circ$.

Hence, suppose that $a\in M_\circ$ and $h\in(T\times G)_\circ$
satisfy the identity $ah=0$; in this case, a simple inspection
shows that either $a=0$ or $h=0$; condition (F1) follows
straightforwardly. To check (F2), say
that $a_1\cdot(s_1,g)=a_2\cdot(s_2,g_2)\neq 0$, for some
$a_i\in M$, $s_i\in S$, $g_i\in G$ ($i=1,2$); since (F2)
holds for $S$, we may find $b_1,b_2\in M$ and $t\in S$ such
that $a_1b_1=a_2b_2$ and $b_it=s_i$ ($i=1,2$). Notice that
$g:=\phi(b_1)^{-1}g_1=\phi(b_2)^{-1}g_2$ in $G$, therefore
$b_i\cdot(t,g)=(s_i,g_i)$ for $i=1,2$, whence the contention.

Conversely, suppose that $S\otimes_{1_\circ}G_\circ$ is
$M_\circ$-flat; we wish to show that (F1) and (F2) hold for $S$.
However, suppose that $as=0$ for some $a\in M_\circ$ and $s\in S$
with $s\neq 0$; then $a\cdot(s\otimes e)=0$ (where $e\in G$ is
the neutral element); since (F1) holds for
$S\otimes_{1_\circ}G_\circ$, we deduce that there exist
$b\in M_\circ$ and $t\otimes g\in S\otimes_{1_\circ}G_\circ$ such
that $ba=0$ and $s\otimes e=b\cdot(t\otimes g)=bt\otimes\phi(b)g$;
this implies that $bt=s$, so (F1) holds for $S$, as sought.

Lastly, suppose that $a_1s_1=a_2s_2\neq 0$ for some $a_i\in M$
and $s_i\in S$ ($i=1,2$). It follows that
$a_1\cdot(s_1\otimes\phi(a_2))=a_2\cdot(s_2\otimes\phi(a_1))$.
By applying condition (F2) to this identity in
$S\otimes_{1_\circ}G_\circ$, we deduce that the same condition
holds also for $S$.
\end{pfclaim}

Next, we observe that there is a natural isomorphism
of $N_\circ$-modules :
\set\begin{equation}\label{eq_iso-times}
N_\circ\otimes_{M_\circ}(S\otimes_{1_\circ}N^\gp_\circ)\isom
(N_\circ\otimes_{M_\circ}S)\otimes_{1_\circ}N_\circ^\gp
\qquad
n\otimes(s\otimes g)\mapsto(n\otimes s)\otimes\phi(n)g
\end{equation}
whose inverse is given by the rule :
$(n\otimes s)\otimes g\mapsto n\otimes(s\otimes\phi(n)^{-1}g)$
for every $n\in N$, $s\in S$, $g\in N^\gp_\circ$.
We leave to the reader the verification that these maps
are well defined, and they are inverse to each other.
In view of \eqref{eq_iso-times} and claim \ref{cl_prod_with-G},
we may then replace $S$ by $S\otimes_{1_\circ}N^\gp_\circ$, which
allows to assume that $S$ is an integral pointed $M_\circ$-module
and $N_\circ\otimes_{M_\circ}S$ is an integral pointed
$N_\circ$-module. In this case, in view of proposition
\ref{prop_Lazard} we know that
$\Z\La N_\circ\otimes_{M_\circ}S\Ra$ is a flat
$\Z\La N_\circ\Ra$-module, and it suffices to show that
$\Z\La S\Ra$ is a flat $\Z\La M_\circ\Ra$-module.

However, under our assumptions, the ring homomorphism
$\Z\La M_\circ\Ra\to\Z\La N_\circ\Ra$ is finite and injective,
so the assertion follows from \cite[Part II, Th.1.2.4]{Gr-Ra}
and \eqref{eq_yet-more-nonsense}.
\end{proof}

\begin{lemma}\label{lem_exc-satura}
Let $M$ be an integral monoid, and $S\subset M$ a submonoid. We have :
\begin{enumerate}
\item
The natural map $S^{-1}M^\sat\to(S^{-1}M)^\sat$ is an isomorphism.
\item
If $M$ is saturated, then $M/S$ is saturated, and if $S$
is a group, also the converse holds.
\item
The inclusion map $M\to M^\sat$ is a local morphism.
\item
The inclusion $M\subset M^\sat$ induces a natural bijection :
$$
\Spec\,M^\sat\isom\Spec\,M.
$$
\end{enumerate}
\end{lemma}
\begin{proof} (i) and (iii) are left to the reader. For (ii),
the natural isomorphism $S^{-1}M/S^\gp\isom M/S$, together
with (i), reduces to the case where $S$ is a group, in which
case it suffices to remark that $(M/S)^\sat=M^\sat/S$.
Lastly, to show (iv) it suffices to prove that, for any face
$F\subset M$, the submonoid $F^\sat\subset M^\sat$ is a face,
and $F^\sat\cap M=F$, and that every face of $M^\sat$ is of
this form. These assertions are easy exercises, which
we leave as well to the reader.
\end{proof}

\begin{lemma}\label{lem_decomp-sats}
Let $M$ be a saturated monoid such that $M^\sharp$ is fine.
Then there exists an isomorphism of monoids :
$$
M\isom M^\sharp\times M^\times
$$
and if $M$ is fine, $M^\times$ is a finitely generated abelian
group. Moreover, the projection $M\to M^\sharp$ induces a bijection :
$$
\Spec\,M^\sharp\isom\Spec\,M
\quad : \quad
\fp\mapsto\fp\times M^\times.
$$
\end{lemma}
\begin{proof} Under the stated assumptions, $G:=M^\gp/M^\times$ is
a free abelian group of finite rank, hence the projection $M^\gp\to
G$ admits a splitting $\sigma:G\to M^\gp$. Set $M_0:=M\cap\sigma(G)$;
it is easily seen that $M=M_0\times M^\times$, whence
$M_0\simeq M^\sharp$. If $M$ is fine, $M^\gp$ is a finitely generated
abelian group, hence the same holds for its direct factor $M^\times$.
The last assertion can be proven directly, or can be regarded as a
special case of lemma \ref{lem_spec-quots}(i.b).
\end{proof}

\begin{definition}\label{def_satur-morph}
Let $\phi:M\to N$ be a morphism of integral monoids. 
\begin{enumerate}
\item
We say that $\phi$ is {\em $k$-saturated} (for some integer $k>0$),
if the push-out $P\otimes_MN$ is integral and $k$-saturated, for
every morphism $M\to P$ with $P$ integral and $k$-saturated.
\item
We say that $\phi$ is {\em saturated}, if the following holds.
For every morphism of monoids $M\to P$ such that $P$ is integral
and saturated, the monoid $P\otimes_MN$ is also integral and
saturated.
\end{enumerate}
\end{definition}

Clearly, if $\phi$ is $k$-saturated for every integer $k>0$,
then $\phi$ is saturated.

\begin{lemma}\label{lem_little}
Let $\phi:M\to N$ be a morphism of integral monoids,
and $S\subset M$, $T\subset N$ two submonoids such
that $\phi(S)\subset T$.  The following holds :
\begin{enumerate}
\item
The localization map $M\to S^{-1}M$ is saturated.
\item
If $\phi$ is saturated, the same holds for the morphisms
$S^{-1}M\to T^{-1}N$ and $M/S\to N/T$ induced by $\phi$.
\item
If $S$ and $T$ are two groups, then $\phi$ is saturated if
and only if the same holds for the induced morphism $M/S\to N/T$.
\item
If $\phi$ is saturated, then the natural map
$N\otimes_MM^\sat\to N^\sat$ is an isomorphism.
\end{enumerate}
\end{lemma}
\begin{proof}(i) follows from the standard isomorphism :
$S^{-1}M\otimes_MN\isom\phi(S)^{-1}N$, together with lemma
\ref{lem_exc-satura}(i). Next, let $M/S\to P$ and
$S^{-1}M\to Q$ be morphisms of monoids. Then
$$
P\otimes_{M/S}N/T\simeq(P\otimes_MN)/T
\qquad\text{and}\qquad
Q\otimes_{S^{-1}M}T^{-1}N\simeq T^{-1}(Q\otimes_MN)
$$
(lemma \ref{lem_localize}) so assertions (ii) and (iii)
follow from lemma \ref{lem_exc-satura}(i,ii).

(iv) follows by comparing the universal properties.
\end{proof}

\begin{example}\label{ex_Rees-satura}
(i)\ \
Let $M$ be an integral monoid, $I\subset M$ an ideal,
and consider again the Rees monoid $\sR(M,I)$ of the
proof of theorem \ref{th_flat-crit-for-mnds}. Clearly
$\sR(M,I)$ is an integral monoid. However, easy examples
show that, even when $M$ is saturated, $\sR(M,I)$ is not
generally saturated. More precisely, the following holds.
For every ideal $J\subset M$, set 
$$
J^\sat:=\{a\in M^\gp~|~a^n\in J^n\text{ for some integer $n>0$}\}
$$
where $J^n$ denotes the $n$-th power of $J$ in the
monoid $(\cP,\cdot)$ of \eqref{subsec_toric}. Then $J^\sat$
is an ideal of $M^\sat$. With this notation, we have the identity :
$$
\sR(M,I)^\sat=\bigoplus_{n\in\N}(I^n)^\sat.
$$
(Verification left to the reader.)

(ii)\ \
For instance, take $M:=\Q^{\oplus 2}_+$, and let $I\subset M$
be the ideal consisting of all pairs
$(x,y)$ such that $x+y>1$. Then
$I^n=\{(x,y)\in\Q_+^{\oplus 2}~|~x+y>n\}$, and it is easily
seen that $I^n=(I^n)^\sat$ for every $n\in\N$, hence in this
case $\sR(M,I)$ is saturated, in view of (ii).
It is easily seen that $\sR(M,I)$ does not fulfill condition
(F2) of lemma \ref{lem_Lazard}, hence the natural inclusion
map $i:M\to\sR(M,I)$ is not flat (theorem \ref{th_always-exact}),
hence it is not an integral morphism, according to theorem
\ref{th_flat-crit-for-mnds}. On the other hand, we have :
\end{example}

\begin{lemma}\label{lem_contras}
With the notation of example {\em\ref{ex_Rees-satura}(ii)}, the
morphism $i$ is saturated.
\end{lemma}
\begin{proof} We prefer to work with the multiplicative
notation, so we shall argue with the monoid
$(\exp\Q^{\oplus 2}_+,\cdot)$ (see \eqref{sec_toric}).
Indeed, let $\phi:M\to P$ be any morphism of monoids, with
$P$ saturated. Clearly $\sR(M,I)$ is the direct sum of the
$P$-modules $P\otimes_MI^n$, for all $n\in\N$.

\begin{claim}\label{cl_counter-satura}
The natural map $P\otimes_MI^n\to I^nP$ is an isomorphism,
for every $n\in\N$
\end{claim}
\begin{pfclaim} Indeed, the map is obviously surjective.
Hence, suppose that $a_1x_1=a_2x_2$ for some $a_1,a_2\in I^n$,
$x_1,x_2\in P$ such that $x_1\otimes a_1=x_2\otimes a_2$.
For every $\theta\in\Q$ with $0\leq\theta\leq 1$, set
$a_\theta:=a_1^\theta\cdot a_2^{1-\theta}$ and notice that
$$
x_\theta:=a_1a_\theta^{-1} x_1=a_2a_\theta^{-1}x_2\in P^\gp
$$
and if $N\in\N$ is large enough, so that $N\theta\in\N$,
then $x_\theta^N=x_1^{N\theta}x_2^{N\theta}\in P$, hence
$x_\theta\in P$.

Next, for any $(a,b),(a',b')\in\Q^{\oplus 2}_+$, set
$(a,b)\vee(a',b'):=(\min(a,a'),\min(b,b'))$. Choose
an increasing sequence $0:=\theta_0<\theta_1<\cdots<\theta_n:=1$
of rational numbers, such that
$$
b_i:=a_{\theta_i}\vee a_{\theta_{i+1}}\in I^n
\qquad
\text{for every $i=0,\dots,n-1$}.
$$
Then there exist $c_i,d_i\in\Q^{\oplus 2}_+$ such that
$a_{\theta_i}=b_ic_i$ and $a_{\theta_{i+1}}=b_id_i$ for
$i=0,\dots,n-1$. We may then compute in $I^n\otimes_MP$ :
$$
a_{\theta_i}\otimes x_{\theta_i}=b_ic_i\otimes x_{\theta_i}
=b_i\otimes c_i x_{\theta_i}
\quad\text{and likewise}\quad
a_{\theta_{i+1}}\otimes x_{\theta_{i+1}}=
b_i\otimes d_i x_{\theta_{i+1}}.
$$
By construction, we have $c_i x_{\theta_i}=d_i x_{\theta_{i+1}}$
for $i=0,\dots,n-1$, whence the contention.
\end{pfclaim}

In view of claim \ref{cl_counter-satura}, we are reduced to
showing that $\sR(P,IP)$ is saturated, and by example
\ref{ex_Rees-satura}, this comes down to proving that
$(I^nP)^\sat=I^nP$ for every $n\in\N$. However, say that
$x\in P^\gp$, and $x^n=a_1x_1\cdots a_nx_n$ for some $a_i\in I$
and $x_i\in P$; set $a:=(a_1\cdots a_n)^{1/n}$, and notice
that $a\in I$. Then $x^n=a^n\cdot x_1\cdots x_n$, so that
$xa^{-1}\in P$, and finally $x\in IP$, as required.
\end{proof}

Lemma \ref{lem_contras} shows that a saturated morphism is not
necessarily integral. Notwithstanding, we shall see later that
integrality holds for an important class of saturated morphisms
(corollary \ref{cor_satu-gives-flat}). Now we wish to globalize
the class of saturated morphisms, to an arbitrary topos.
Of course, we could define the notion of saturated morphism
of $T$-monoids, just by repeating word by word definition
\ref{def_satur-morph}(ii). However, it is not clear that the
resulting condition would be of a type which can be checked
on stalks, in the sense of definition \ref{def_T-point}(v).
For this reason, we prefer to proceed as in Tsuji's work \cite{Tsu}.

\begin{lemma}\label{lem_gener-exact}
Let $T$ be a topos, $\phi:\underline M\to\underline N$ and
$\psi:\underline N\to\underline P$ two morphisms of integral
$T$-monoids. We have:
\begin{enumerate}
\item
If $\phi$ and $\psi$ are exact, the same holds for $\psi\circ\phi$.
\item
If $\psi\circ\phi$ is exact, the same holds for $\phi$.
\item
Consider a commutative diagram of integral $T$-monoids :
\set\begin{equation}\label{eq_diag-mons}
{\diagram
\underline M \ar[r]^-\phi \ar[d]_\psi & \underline N \ar[d]^{\psi'} \\
\underline M' \ar[r]^-{\phi'} & \underline N'.
\enddiagram}
\end{equation}
Then the following holds :
\begin{enumerate}
\item
If \eqref{eq_diag-mons} is a cartesian diagram and $\phi'$ is exact,
then $\phi$ is exact.
\item
If \eqref{eq_diag-mons} is cocartesian (in the category
$\IntMnd_T$) and $\phi$ is exact, then $\phi'$ is exact.
\end{enumerate}
\end{enumerate}
\end{lemma}
\begin{proof} For all these assertions, remark \ref{rem_satura}(ii)
easily reduces to the case where $T=\Set$, which therefore we assume
from start. Now, (i) and (ii) are left to the reader.

(iii.a): Let $x\in M^\gp$ such that $\phi^\gp(x)\in N$; hence
$(\phi'\circ\psi)^\gp(x)=\psi'(\phi^\gp(x))\in N'$, and therefore
$\psi^\gp(x)\in M'$, since $\phi'$ is exact.
It follows that $x\in M$, so $\phi$ is exact.

(iii.b): Let $x\in (M')^\gp$ such that $(\phi')^\gp(x)\in N'$; then
we may write $(\phi')^\gp(x)=\phi'(y)\cdot\psi'(z)$ for some
$y\in M'$ and $z\in N$. Therefore, $(\phi')^\gp(xy^{-1})=\psi'(z)$;
since the functor $P\mapsto P^\gp$ commutes with colimits, it follows
that we may find $w\in M^\gp$ such that $\psi^\gp(w)=xy^{-1}$ and
$\phi^\gp(w)=z$. Since $\phi$ is exact, we deduce that $w\in M$,
therefore $xy^{-1}\in M'$, and finally $x\in M'$, whence the contention.
\end{proof}

\sset\subsubsection{}\label{subsec_These-consider}
Let $T$ be a topos; for two morphisms
$\underline P\leftarrow\underline M\to\underline N$ of integral
$T$-monoids, we set
$$
\underline N\otimesint_M\underline P:=
(\underline N\otimes_M\underline P)^\intg
$$
which is the push-out of these morphisms, in the category
$\IntMnd_T$. Notice that, for every morphism $f:T'\to T$
of topoi, the natural morphism of $T'$-monoids
\set\begin{equation}\label{eq_int-push-out}
f^*(\underline N\otimesint_M\underline P)\to
f^*\underline N\otimesint_{f^*M}f^*\underline P
\end{equation}
is an isomorphism (by lemmata \ref{lem_T-satura}(i) and
\ref{lem_concern}(i)).

Let now $\phi:\underline M\to\underline N$ be a morphism of
integral $T$-monoids.
For any integer $k>0$, let $\bek_M$ and $\bek_N$ be the
$k$-Frobenius maps of $M$ and $N$ (definition
\ref{def_exact-phi}(ii)), and consider the cocartesian
diagram :
\set\begin{equation}\label{eq_motivate-sat}
{\diagram \underline M \ar[rr]^-{\bek_M} \ar[d]_\phi & & 
          \underline M \ar[d]^{\phi'} \\
          \underline N \ar[rr]^-{\bek_M\otimesint_M\one_N} & &
\underline P.
\enddiagram}
\end{equation}
The endomorphism $\bek_N$ factors through $\bek_M\otimesint_M\one_N$
and a unique morphism $\beta:\underline P\to\underline N$ such that
$\beta\circ\phi'=\phi$. A simple inspection shows that
$\bek_P=(\bek_M\otimes_M\one_N)\circ\beta$.
Now, if $\bek_M$ is exact, then the
same holds for $\bek_M\otimesint_M\one_N$, in view of lemma
\ref{lem_gener-exact}(iii.b). Hence, if $\underline P$ is
$k$-saturated, then $\beta$ is exact (lemma
\ref{lem_gener-exact}(ii)), and if $\underline M$ is $k$-saturated,
also the converse holds. Likewise, if $\underline M$ is
$k$-saturated and $\beta$ is exact, then $\underline N$ is
$k$-quasi-saturated.

These considerations motivate the following :

\begin{definition}
Let $T$ be a topos, $\phi:\underline M\to\underline N$ a morphism
of integral $T$-monoids.
\begin{enumerate}
\item
A commutative diagram of integral $T$-monoids :
$$
\xymatrix{ \underline M \ar[r]^-\phi \ar[d] & \underline N \ar[d] \\
           \underline M' \ar[r]^-{\phi'} & \underline N'
}$$
is called an {\em exact square}, if the induced morphism
$\underline M'\otimesint_M\underline N\to\underline N'$ is exact.
\item
$\phi$ is said to be {\em $k$-quasi-saturated\/} if the
commutative diagram :
\set\begin{equation}\label{eq_n-quasi-sat}
{\diagram
\underline M \ar[r]^-\phi \ar[d]_{\bek_M} & 
\underline N \ar[d]^{\bek_N} \\
\underline M \ar[r]^-\phi & \underline N
\enddiagram}
\end{equation}
is an exact square (the vertical arrows are the $k$-Frobenius maps).
\item
$\phi$ is said to be {\em quasi-saturated\/} if it is
$k$-quasi-saturated for every integer $k>0$.
\end{enumerate}
\end{definition}

\begin{proposition}\label{prop_ex-morphis}
Let $T$ be a topos.
\begin{enumerate}
\item
If \eqref{eq_diag-mons} is an exact square, and
$\underline M\to\underline P$ is any morphism of integral
$T$-monoids, the square
$\underline P\otimesint_M\eqref{eq_diag-mons}$ is exact.
\item
Consider a commutative diagram of integral $T$-monoids :
\set\begin{equation}\label{eq_two-squares}
{\diagram
\underline M \ar[r]^-{\phi_1} \ar[d]_\psi &
\underline N \ar[d] \ar[r]^-{\phi_2} &
\underline P \ar[d]^{\psi'} \\
\underline M' \ar[r]^-{\phi'_1} &
\underline N' \ar[r]^-{\phi'_2} & \underline P'
\enddiagram}
\end{equation}
and suppose that the left and right square subdiagrams of
\eqref{eq_two-squares} are exact. Then the same holds for
the square diagram :
$$
\xymatrix{
\underline M \ar[rr]^-{\phi_2\circ\phi_1} \ar[d]_\psi & &
\underline P \ar[d]^{\psi'} \\
\underline M' \ar[rr]^-{\phi_2'\circ\phi_1'} & &
\underline P'.
}$$
\end{enumerate}
\end{proposition}
\begin{proof} (i): We have a commutative diagram :
$$
\xymatrix{
(\underline P\otimesint_M\underline M')\otimesint_M\underline N 
\ar[r]^-\sigma \ar[d]_\omega &
\underline P\otimesint_M(\underline M'\otimesint_M\underline N)
\ar[d]^{\one_P\otimesint_M\alpha} \\
(\underline P\otimesint_M\underline M')\otimesint_P
(\underline P\otimesint_M\underline N)
\ar[r]^-\beta & \underline P\otimesint_M\underline N'
}$$
where $\omega$, $\sigma$ are the natural isomorphisms, and $\beta$
and $\alpha:\underline M'\otimesint_M\underline N\to\underline N'$
are the natural maps.
By assumption, $\alpha$ is exact, hence the same holds for
$\one_P\otimesint_M\alpha$, by lemma \ref{lem_gener-exact}(iii.b).
So $\beta$ is exact, which is the claim.

(ii): Let
$\alpha:\underline M'\otimesint_M\underline N\to\underline  N'$ and
$\beta:\underline N'\otimesint_N\underline P\to\underline P'$ be the
natural maps; by assumptions, these are exact morphisms. However, we
have a natural commutative diagram :
$$
\xymatrix{
(\underline M'\otimesint_M\underline N)\otimesint_N\underline P
\ar[r]^-\omega \ar[d]_{\alpha\otimesint_N\one_P} &
\underline M'\otimesint_M\underline P \ar[d]^\gamma \\
\underline N'\otimesint_N\underline P \ar[r]^-\beta &
\underline P'
}$$
where $\omega$ is the natural isomorphism, and $\gamma$ is the
map deduced from $\psi'$ and $\phi_2'\circ\phi_1'$. Then the
assertion follows from lemma \ref{lem_gener-exact}(i,iii.b).
\end{proof}

\begin{corollary}\label{cor_exact}
Let $T$ be a topos, $\phi:\underline M\to\underline N$,
$\psi:\underline N\to\underline P$ two morphisms of integral
$T$-monoids, and $h,k>0$ any two integers. The following holds :
\begin{enumerate}
\item
If $\phi$ is both $h$-quasi-saturated and $k$-quasi-saturated, then
$\phi$ is $hk$-quasi-saturated.
\item
If $\phi$ and $\psi$ are $k$-quasi-saturated, the same holds for
$\psi\circ\phi$.
\item
If $\underline M\to\underline P$ is any map of integral $T$-monoids,
and $\phi$ is $k$-quasi-saturated (resp. quasi-saturated), then
the same holds for
$\phi\otimesint_M\one_P:
\underline P\to\underline N\otimesint_M\underline P$.
\item
Let $\underline S\subset\phi^{-1}(\underline N^\times)$ be a
$T$-submonoid. Then $\phi$ is $k$-quasi-saturated (resp.
quasi-saturated) if and only if the same holds for the induced
map $\phi_S:\underline S^{-1}\underline M\to\underline N$.
\item
Let $\underline G\subset\Ker\,\phi$ be a subgroup. Then $\phi$
is $k$-quasi-saturated (resp. quasi-saturated) if and only
if the same holds for the induced map
$\bar\phi:\underline M/\underline G\to\underline N$.
\item
$\phi$ is quasi-saturated if and only if it is $p$-quasi-saturated
for every prime number $p$.
\item
$\underline M$ is $k$-saturated (resp. saturated) if and only
if the unique morphism of $T$-monoids $\{1\}\to\underline M$ is
$k$-quasi-saturated (resp. quasi-saturated).
\item
If $\phi$ is $k$-quasi-saturated (resp. quasi-saturated), and
$\underline M$ is $k$-saturated (resp. saturated), then
$\underline N$ is $k$-saturated (resp. saturated).
\item
If $\underline M$ is integral, and $\underline N$ is a $T$-group,
$\phi$ is quasi-saturated.
\end{enumerate}
\end{corollary}
\begin{proof} (i) and (ii) are straightforward consequences
of proposition \ref{prop_ex-morphis}(ii).

To show (iii), set
$\underline P':=\underline N\otimesint_M\underline P$ and let
us remark that we have a commutative diagram
$$
\xymatrix{
\underline P \ar[d]_-{\phi\otimesint_M\one_P} \ar[r] &
\underline P_1 \ar[r] \ar[d] &
\underline P \ar[d]^{\phi\otimesint_M\one_P} \\
\underline P' \ar[r] & \underline P_2 \ar[r] &
\underline P'
}$$
such that :
\begin{itemize}
\item
the composition of the top (resp. bottom) arrows is the
$k$-Frobenius map
\item
the left square subdiagram is
$\eqref{eq_n-quasi-sat}\otimesint_M\underline P$
\item
the right square subdiagram is cocartesian (hence exact).
\end{itemize}
then the assertion follows from proposition \ref{prop_ex-morphis}.

(iv): Suppose that $\phi$ is $k$-quasi-saturated. Then the
same holds for $\underline S^{-1}\phi:\underline S^{-1}M\to
\underline S^{-1}M\otimesint_M\underline N$, according to (iii).
However, we have a natural isomorphism
$\underline S^{-1}\underline M\otimesint_M\underline N\isom
\phi(\underline S)^{-1}\underline N=N$, so $\phi_S$ is
$k$-quasi-saturated.

Conversely, suppose that $\phi_S$ is $k$-quasi-saturated.
By (ii), in order to prove that the same holds for $\phi$,
it suffices to show that the localisation map
$\underline M\to\underline S^{-1}\underline M$ is
$k$-quasi-saturated. But this is clear, since
\eqref{eq_n-quasi-sat} becomes cocartesian if we take for
$\phi$ the localisation map.

(v): To begin with, $\underline M/\underline G$ is an integral
monoid, by lemma \ref{lem_integral-quot}. Suppose that $\phi$
is $k$-quasi-saturated. Arguing as in the proof of (iv) we see
that the natural map
$\underline N\to(\underline M/\underline G)\otimesint_M\underline N$
is an isomorphism, hence $\bar\phi$ is $k$-quasi-saturated, by (iii).

Conversely, suppose that $\bar\phi$ is $k$-quasi-saturated.
By (ii), in order to prove that $\phi$ is saturated, it suffices
to show that the same holds for the projection
$\underline M\to\underline M/\underline G$.
By (iii), we are further reduced to showing that the unique map
$\underline G\to\{1\}$ is $k$-quasi-saturated, which is trivial.

(vi) is a straightforward consequence of (i). Assertion (vii) can be
verified easily on the definitions. Next, suppose that $\phi$ is
quasi-saturated and $\underline M$ is saturated. By (ii) and (vii)
it follows that the unique morphism $\{1\}\to\underline N$ is
quasi-saturated, hence (vii) implies that $\underline N$ is saturated,
so (viii) holds.

(ix): In view of (iv), we are reduced to the case where $\underline M$
is also a group, in which case the assertion is obvious.
\end{proof}

\begin{proposition}\label{prop_motivate-def}
Let $\phi:M\to N$ be an integral morphism of integral monoids,
$k>0$ an integer. The following conditions are equivalent :
\begin{enumerate}
\alphaenu
\item
$\phi$ is $k$-quasi-saturated.
\item
$\phi$ is $k$-saturated.
\item
The push-out $P$ of the cocartesian diagram \eqref{eq_motivate-sat},
is $k$-saturated.
\end{enumerate}
\end{proposition}
\begin{proof}(a)$\Rightarrow$(b) by virtue of corollary
\ref{cor_exact}(iii,viii), and trivially (b)$\Rightarrow$(c).
Lastly, the implication (c)$\Rightarrow$(a) was already
remarked in \eqref{subsec_These-consider}.
\end{proof}

\begin{corollary}
Let $\phi:M\to N$ be an integral morphism of integral monoids.
Then $\phi$ is quasi-saturated if and only if it is saturated.
\qed\end{corollary}

Proposition \ref{prop_motivate-def} motivates the following :

\begin{definition}\label{def_satur-top}
Let $T$ be a topos, $\phi:\underline M\to\underline N$ a
morphism of integral $T$-monoids, $k>0$ an integer. We say
that $\phi$ is {\em $k$-saturated} (resp. {\em saturated})
if $\phi$ is integral and $k$-quasi-saturated (resp. and
quasi-saturated).
\end{definition}

In case $T=\Set$, we see that an integral morphism of
integral monoids is saturated in the sense of definition
\ref{def_satur-top}, if and only if it is saturated in the
sense of the previous definition \ref{def_satur-morph}.
We may now state :

\begin{corollary}\label{cor_satur-check}
Let\/ $\bP(T,\phi)$ be the property ``$\phi$ is an integral
(resp. exact, resp. $k$-saturated, resp. saturated) morphism
of integral $T$-monoids'' (for a topos $T$). Then $\bP$ can
be checked on stalks. (See definition {\em\ref{def_T-point}(v)}.)
\end{corollary}
\begin{proof} The fact that integrality of a $T$-monoid can
be checked on stalks, has already been established in lemma
\ref{lem_concern}(ii). For the property ``$\phi$ is an integral
morphism'' (between integral $T$-monoids), it suffices to apply
theorem \ref{th_flat-crit-for-mnds} and proposition
\ref{prop_finally-flat}.

For the property ``$\phi$ is exact'' one applies lemma
\ref{lem_concern}(iii). From this, and from \eqref{eq_int-push-out}
one deduces that also the properties of $k$-saturation
and saturation can be checked on stalks.
\end{proof}

\begin{lemma}\label{lem_persist-integr}
Let $f:M\to N$ be a morphism of integral monoids. We have :
\begin{enumerate}
\item
If $f$ is exact, then $f$ is local.
\item
Conversely, if $f$ is an integral and local morphism, then
$f$ is exact.
\end{enumerate}
\end{lemma}
\begin{proof}(i) is left to the reader as an exercise.

(ii): Suppose $x\in M^\gp$ is an element such that
$b:=f_S(x)\in N$. Write $x=z^{-1}y$ for certain $y,z\in M$;
therefore $b\cdot f(z)=f(y)$ holds in $N$, and theorems
\ref{th_flat-crit-for-mnds}, \ref{th_always-exact} imply
that there exist $c\in N$ and $a_1,a_2\in M$, such that
$1=c\cdot f(a_1)$, $b=c\cdot f(a_2)$ and $ya_1=za_2$.
Since $f$ is local, we deduce that $a_1\in M^\times$,
hence $x=a_1^{-1}a_2$ lies in $M$.
\end{proof}

\begin{proposition}\label{prop_crit-saturated}
Let $f:M\to N$ be an integral map of integral monoids, $k>0$ an
integer,  and $N=\bigoplus_{\gamma\in\Gamma}N_\gamma$ the
$f$-grading of $N$ (see remark {\em\ref{rem_why-not}(iii)}).
Then the following conditions are equivalent :
\begin{enumerate}
\alphaenu
\item
$f$ is $k$-quasi-saturated.
\item
$N_{k\gamma}=N_\gamma^k$\ \ for every $\gamma\in\Gamma$. (Here the
$k$-th power of subsets of $N$ is taken in the monoid $\cP(N)$,
as in \eqref{subsec_toric}.)
\end{enumerate}
\end{proposition}
\begin{proof} (a)$\Rightarrow$(b): Let $\pi:N^\gp\to\Gamma$ be
the projection, suppose that $y\in N_{k\gamma}$ for some
$\gamma\in\Gamma$, and pick $x\in N^\gp$ such that $\pi(x)=\gamma$.
This means that $y=x^k\cdot f^\gp(z)$ for some $z\in M^\gp$.
By (a), it follows that we may find a pair $(a,b)\in M\times N$
and an element $w\in M^\gp$, such that $(aw^{-k},bw)=(z,x)$ in
$M^\gp\times N^\gp$. Especially, $b,b\cdot f(a)\in N_\gamma$,
and consequently $y=b^{k-1}\cdot(b\cdot f(a))\in N_{k\gamma}$,
as stated.

(b)$\Rightarrow$(a): Let $S:=f^{-1}(N^\times)$; by lemmata
\ref{lem_persist-integr}(ii) and \ref{lem_obvious-int}(iii), the induced
map $f_S:S^{-1}M\to N$ is exact and integral. Moreover, by corollary
\ref{cor_exact}(iv), $f$ is $k$-quasi-saturated if and only if the
same holds for $f_S$, and clearly the $f$-grading of $N$ agrees with
the $f_S$-grading. Hence we may replace $f$ by $f_S$ and assume from
start that $f$ is exact. In this case, $G:=\Ker\,f^{\gp}\subset M$,
and corollary \ref{cor_exact}(v) says that $f$ is $k$-quasi-saturated
if and only if the same holds for the induced map $\bar f:M/G\to N$;
moreover $\bar f$ is still integral, since it is deduced from $f$ by
push-out along the map $M\to M/G$. Hence we may replace $f$ by
$\bar f$, thereby reducing to the case where $f$ is injective.
Also, $f$ is flat, by theorem \ref{th_flat-crit-for-mnds}.
The assertion boils down to the following. Suppose that
$(x,y)\in M^\gp\times N^\gp$ is a pair such that
$a:=f^\gp(x)\cdot y^k\in N$; we have to show that there exists
a pair $(m,n)\in M\times N$ whose class in the push-out $P$ of
the diagram
$$
N\xleftarrow{\ f\ }M\xrightarrow{\ \bek_M\ }M
$$
agrees with the image of $(x,y)$ in $P^\gp$. However, set
$\gamma:=\pi(y)$, and notice that $a\in N_{k\gamma}$, hence
we may write $a=n_1\cdots n_k$ for certain
$n_1,\dots,n_k\in N_\gamma$. Then, according to remark
\ref{rem_why-not}(iv), we may find $n\in N_\gamma$ and
$x_1,\dots,x_k\in M$ such that $n_i=n\cdot f(x_i)$ for every
$i=1,\dots,k$. It follows that $y\cdot n^{-1}\in f^\gp(M^\gp)$, say
$y=n\cdot f(z)$ for some $z\in M^\gp$. Set $m:=x_1\cdots x_k$; then
$(x,y)=(x,n\cdot f(z))$ represents the same class as
$(x\cdot z^k,n)$ in $P^\gp$. Especially, $f(x\cdot z^k)\cdot
n^k=a=f(m)\cdot n^k$, hence $m=x\cdot n^k$, since $f$ is injective.
The claim follows.
\end{proof}

\begin{corollary}\label{cor_persist-integr}
Let $f:M\to N$ be an integral and local morphism of integral
and sharp monoids. Then :
\begin{enumerate}
\item
$f$ is exact and injective.
\item
If furthermore, $f$ is saturated, then $\Coker\,f^\gp$ is a
torsion-free abelian group.
\end{enumerate}
\end{corollary}
\begin{proof}(i): It follows from lemma \ref{prop_crit-saturated}
that $f$ is exact. Next, suppose that $f(a_1)=f(a_2)$ for some
$a_1,a_2\in M$.
By theorem \ref{th_flat-crit-for-mnds} and \ref{th_always-exact},
it follows that there exist $b_1,b_2\in M$ and $t\in N$ such
that $1=f(b_1)t=f(b_2)t$ and $a_1b_1=a_2b_2$. Since $N$
is sharp, we deduce that $f(b_1)=f(b_2)=1$, and since $f$
is local and $M$ is sharp, we get $b_1=b_2=1$, thus $a_1=a_2$,
whence $x=1$, which is the contention.

(ii): Let us endow $N$ with its $f$-grading. Suppose now that
$g\in G$ is a torsion element, and say that $g^k=1$ in $G$ for
some integer $k>0$; by propositions \ref{prop_crit-saturated}
and \ref{prop_motivate-def}, we then have $N_1=N_g^k$.
Especially, there exist $a_1,\dots,a_k\in N_g$,
such that $a_1\cdots a_k=1$. Since $N$ is sharp, we must
have $a_i=1$ for $i=1,\dots,k$, hence $g=1$.
\end{proof}

\begin{corollary}\label{coro_satur-face}
Let $\phi:M\to N$ be a morphism of integral monoids, $F\subset N$
any face, and $\phi_F:\phi^{-1}F\to F$ the restriction of $\phi$.
\begin{enumerate}
\item
If $\phi$ is flat, then the same holds for $\phi_F$, and
the induced map $\Coker\,\phi_F^\gp\to\Coker\,\phi^\gp$ is
injective.
\item
If $\phi$ is saturated, the same holds for $\phi_F$.
\end{enumerate}
\end{corollary}
\begin{proof}(i): The fact that $\phi_F$ is flat, is a special
case of corollary \ref{cor_flat-face}. The assertion about
cokernels boils down to showing that the induced diagram of
abelian groups :
$$
\xymatrix{
(\phi^{-1}F)^\gp \ar[r]^-{\phi^\gp_F} \ar[d] & F^\gp \ar[d] \\
M^\gp \ar[r]^-{\phi^\gp} & N^\gp
}$$
is cartesian. However, say that $\phi^\gp(a_1a_2^{-1})=f_1f_2^{-1}$
for some $a_1,a_2\in M$ and $f_1,f_2\in F$. This means that
$\phi(a_1)f_2=\phi(a_2)f_1$ in $N$; by condition (F2) of
theorem \ref{th_always-exact}, we deduce that there exist
$b_1,b_2\in M$ and $t\in N$ such that $b_2a_1=b_1a_2$ and
$f_i=t\cdot\phi(b_i)$ ($i=1,2$). Since $F$ is a face, it
follows that $b_1,b_2\in\phi^{-1}F$, hence
$a_1a_2^{-1}=b_1b_2^{-1}\in(\phi^{-1}F)^\gp$, as required.

(ii): By (i), the morphism $\phi_F$ is integral, hence it
suffices to show that $\phi_F$ quasi-saturated (proposition
\ref{prop_motivate-def}), and to this aim we shall apply the
criterion of proposition \ref{prop_crit-saturated}. Indeed,
let $N=\bigoplus_{\gamma\in\Gamma}N_\gamma$ (resp.
$F=\bigoplus_{\gamma\in\Gamma_{\!F}}F_\gamma$) be the
$\phi$-grading (resp. the $\phi_F$-grading); according to (i),
the induced map $j:\Gamma_{\!F}\to\Gamma$ is injective. This means
that $F_\gamma=F\cap N_{j(\gamma)}$ for every $\gamma\in\Gamma_{\!F}$.
Hence, for every integer $k>0$ we may compute
$$
F_{k\gamma}=F\cap N^k_{j(\gamma)}=(F\cap N_{j(\gamma)})^k=
F_\gamma^k
$$
where the first identity follows by applying proposition
\ref{prop_crit-saturated} (and proposition \ref{prop_motivate-def})
to the saturated map $\phi$, and the second identity holds because
$F$ is a face of $N$.
\end{proof}

\subsection{Polyhedral cones}\label{sec_polyhedra}
Fine monoids can be studied by means of certain combinatorial
objects, which we wish to describe. Part of the material that
follows is borrowed from \cite{Ful}. Again,
{\em all the monoids in this section are non-pointed.}

\sset\subsubsection{}\label{subsec_convex-cone}
Quite generally, a {\em convex cone\/} is a pair $(V,\sigma)$,
where $V$ is a finite dimensional $\R$-vector space, and
$\sigma\subset V$ is a non-empty subset such that :
$$
\R_+\cdot\sigma=\sigma=\sigma+\sigma
$$
where the addition is formed in the monoid $(\cP(V),+)$ as in
\eqref{subsec_toric}, and scalar multiplication by the set
$\R_+$ is given by the rule :
$$
\R_+\cdot S:=\{r\cdot s~|~r\in\R_+, s\in S\} \qquad\text{for every
$S\in\cP(V)$}.
$$
A subset $S\subset\sigma$ is called a {\em ray\/} of $\sigma$, if it
is of the form $\R_+\cdot\{s\}$, for some
$s\in\sigma\!\setminus\!\{0\}$.

We say that $(V,\sigma)$ is a {\em closed convex cone\/} if $\sigma$
is closed as a subset of $V$ (of course, $V$ is here regarded as a
topological space via any choice of isomorphism $V\simeq\R^n$). We
denote by $\langle\sigma\rangle\subset V$ the $\R$-vector space
generated by $\sigma$. To a convex cone $(V,\sigma)$ one assigns the
{\em dual cone\/} $(V^\vee,\sigma^\vee)$, where
$V^\vee:=\Hom_\R(V,\R)$, the dual of $V$, and :
$$
\sigma^\vee:= \{u\in V^\vee~|~\text{$u(v)\geq 0$ for every
$v\in\sigma$}\}.
$$
Also, the {\em opposite cone\/} of $\sigma$ is the cone
$$
-\sigma:=\{-v\in V~|~v\in\sigma\}.
$$
Notice that $\sigma^\vee$ and $-\sigma$ are always closed cones.
Notice as well that the restriction of the addition law of $V$
determines a monoid structure $(\sigma,+)$ on the set $\sigma$.
A map of cones
$$
\phi:(W,\sigma_W)\to(V,\sigma_V)
$$
is an $\R$-linear map $\phi:W\to V$ such that
$\phi(\sigma_W)\subset\sigma_V$. Clearly, the restriction of
$\phi$ yields a map of monoids $(\sigma_W,+)\to(\sigma_V,+)$.
If $S\subset V$ is any subset, we set :
$$
S^\bot:=\{u\in V^\vee~|~\text{$u(s)=0$ for every $s\in S$}\}\subset
V^\vee.
$$

\begin{lemma}\label{lem_ddoble-cone}
Let $(V,\sigma)$ be a closed convex cone, Then, under the natural
identification $V\isom(V^\vee)^\vee$, we have
$(\sigma^\vee)^\vee=\sigma$.
\end{lemma}
\begin{proof} This follows from \cite[Ch.II, \S5, n.3, Cor.5]{BouEVT}.
\end{proof}

\sset\subsubsection{}\label{subsec_conv-polyhedr} A {\em convex
polyhedral cone\/} is a cone $(V,\sigma)$ such that $\sigma$ is of
the form :
$$
\sigma:=\{r_1v_1+\cdots+r_sv_s\in V~|~\text{$r_i\geq 0$ for every
$i\leq s$}\}
$$
for a given set of vectors $v_1,\dots,v_s\in V$, called {\em
generators\/} for the cone $\sigma$. One also says that
$\{v_1,\dots,v_s\}$ is a {\em generating set} for $\sigma$, and that
$\R_+\cdot v_1,\dots,\R_+\cdot v_s$ are {\em generating rays\/} for
$\sigma$. We say that $\sigma$ is a {\em simplicial cone}, if it
is generated by a system of linearly independent vectors.

\begin{lemma}
Let $(V,\sigma)$ be a convex polyhedral cone, $S$ a finite
generating set for $\sigma$. Then :
\begin{enumerate}
\item
For every $v\in\sigma$ there is a subset $T\subset S$ consisting of
linearly independent vectors, such that $v$ is contained in the
convex polyhedral cone generated by $T$.
\item
$(V,\sigma)$ is a closed convex cone.
\end{enumerate}
\end{lemma}
\begin{proof} (i): Let $T\subset S$ be a subset such that
$v$ is contained in the cone generated by $T$; up to replacing $T$
by a subset, we may assume that $T$ is minimal, {\em i.e.} no proper
subset of $T$ generates a cone containing $v$. We claim that $T$
consists of linearly independent vectors. Otherwise, we may find a
linear relation of the form $\sum_{w\in T}a_w\cdot w=0$, for certain
$a_w\in\R$, at least one of which is non-zero; we may then assume
that :
\set\begin{equation}\label{eq_assumpt}
a_w>0
\qquad\text{for at least one vector $w\in T$}.
\end{equation}
Say also that $v=\sum_{w\in T}b_w\cdot w$ with $b_w\in\R_+$; by the
minimality assumption on $T$, we must actually have $b_w>0$ for
every $w\in T$. We deduce the identity : $v=\sum_{w\in
T}(b_w-ta_w)\cdot w$ for every $t\in\R$; let $t_0$ be the supremum
of the set of positive real numbers $t$ such that $b_w-ta_w\geq 0$
for every $w\in T$. In view of \eqref{eq_assumpt} we have
$t_0\in\R_+$; moreover $b_w-t_0a_w\geq 0$ for every $w\in T$, and
$b_w-t_0a_w=0$ for at least one vector $w$, which contradicts the
minimality of $T$.

(ii): In view of (i), we are reduced to the case where $\sigma$ is
generated by finitely many linearly independent vectors, and for
such cones the assertion is clear (details left to the reader).
\end{proof}

\sset\subsubsection{}\label{subsec_rather-see-here}
A {\em face\/} of a convex cone $\sigma$ is a subset of the form
$\sigma\cap\Ker\,u$, for some $u\in\sigma^\vee$. The {\em dimension}
(resp. {\em codimension}) of a face $\tau$ of $\sigma$ is the
dimension of the $\R$-vector space $\langle\tau\rangle$ (resp.
of the $\R$-vector space $\langle\sigma\rangle/\langle\tau\rangle$).
A {\em facet\/} of $\sigma$ is a face of codimension one. Notice
that if $\sigma$ is a polyhedral cone, and $\tau$ is a face of
$\sigma$, then $(V,\tau)$ is also a convex polyhedral cone;
indeed if $S\subset V$ is a generating set for $\sigma$, then
$S\cap\tau$ is a generating set for $\tau$.

\begin{lemma}\label{lem_same-faces}
Let $(V,\sigma)$ be a convex polyhedral cone. Then the faces of
$(V,\sigma)$ are the same as the faces of the monoid $(\sigma,+)$.
\end{lemma}
\begin{proof} Clearly we may assume that $\sigma\neq\{0\}$.
Let $F$ be a face of the polyhedral cone $\sigma$, and pick $u\in
\sigma^\vee$ such that $F=\sigma\cap\Ker\,u$. Then
$\sigma\setminus F=\{x\in\sigma~|~u(x)>0\}$, and this is clearly
an ideal of the monoid $(\sigma,+)$; hence $F$ is a face of
$(\sigma,+)$.

Conversely, suppose that $F$ is a face of $(\sigma,+)$. First,
we wish to show that $F$ is a cone in $V$. Indeed, let $f\in F$,
and $r>0$ any real number; we have to prove that $r\cdot f\in F$.
To this aim, it suffices to show that $(r/N)\cdot f\in F$ for some
integer $N>0$, so that, after replacing $r$ by $r/N$ for $N$ large
enough, we may assume that $0<r<1$. In this case,
$(1-r)\cdot f\in \sigma$, and we have $f=r\cdot f+(1-r)\cdot f$,
so that $r\cdot f\in F$, since $F$ is a face of $(\sigma,+)$.

Next, denote by $W\subset V$ the $\R$-vector space spanned by
$F$. Suppose first that $V=W$, and consider any $m\in\sigma$;
then $m=f_1-f_2$ for some $f_1,f_2\in F$, hence $f_1=f_2+m$. Since
$F$ is a face of $(\sigma,+)$, this implies that $m$ lies in $F$, so
$F=\sigma$ in this case, especially $F$ is a face of the convex
polyhedral cone $\sigma$.

So finally, we may assume that $W$ is a proper subspace of
$V$. In this case, let $N:=\sigma+(-F)$. We notice that
$N\neq V$. Indeed, if $m\in\sigma$ and $-m\in N$, we may write
$-m=m'-f$ for some $m'\in\sigma$ and $f\in F$; hence $f=m+m'$,
and therefore $m\in F$; in other words, $N$ avoids the whole of
$-(\sigma\setminus F)$, which is not empty for $\sigma\neq\{0\}$.

Thus, $N$ is a proper convex cone of $V$. Now, let
$u_1,\dots,u_k\in N^\vee$ be a system of generators of the
$\R$-subspace of $V^\vee$ spanned by $N^\vee$, and set
$u:=u_1+\cdots+u_k$. Suppose that $x\in\sigma\cap\Ker\,u$; then
$-x\in\Ker\,u_i$ for every $i=1,\dots,k$, and therefore
$-x\in N^{\vee\vee}=N$ (lemma \ref{lem_ddoble-cone}). Hence, we
may write $-x=m-f$ for some $m\in\sigma$ and $f\in F$, or
equivalently, $f=m+x$, which shows that $x\in F$. Summarizing,
we have proved that $F=\sigma\cap\Ker\,u$, {\em i.e.} $F$ is
a face of the convex cone $\sigma$.
\end{proof}

\begin{proposition}\label{prop_refer-to-Ful}
Let $(V,\sigma)$ be a convex polyhedral cone. The following holds :
\begin{enumerate}
\item
Any intersection of faces of $\sigma$ is still a face of $\sigma$.
\item
If $\tau$ is a face of $\sigma$, and $\gamma$ is a face of
$(V,\tau)$, then $\gamma$ is a face of $\sigma$.
\item
Every proper face of $\sigma$ is the intersection of the facets that
contain it.
\end{enumerate}
\end{proposition}
\begin{proof} (i): Say that $\tau_i=\sigma\cap\Ker\,u_i$, where
$u_1,\dots,u_n\in\sigma^\vee$. Then
$\bigcap_{i=1}^n\tau_i=\sigma\cap\Ker\sum_{i=1}^n u_i$.

(ii): Say that $\tau=\sigma\cap\Ker\,u$ and
$\gamma=\tau\cap\Ker\,v$, where $u\in\sigma^\vee$ and
$v\in\tau^\vee$. Then, for large $r\in\R_+$, the linear form
$v':=v+ru$ is non-negative on any given finite generating set of
$\sigma$, hence it lies in $\sigma^\vee$, and
$\gamma=\sigma\cap\Ker\,v'$.

(iii): To begin with, we prove the following :

\begin{claim}\label{cl_codim-face}
(i)\ \ Every proper face of $\sigma$ is contained in a facet.
\begin{enumerate}
\addenu
\item
Every face of $\sigma$ of codimension $2$ is the intersection of
exactly two facets.
\end{enumerate}
\end{claim}
\begin{pfclaim} We may assume that $\langle\sigma\rangle=V$.
Let $\tau$ be a face of $\sigma$ of codimension at least two, and
denote by $\bar\sigma$ be the image of $\sigma$ in the quotient
$\bar V:=V/\langle\tau\rangle$; clearly $(\bar V,\bar\sigma)$ is
again a polyhedral cone. Moreover, choose $u\in\sigma^\vee$ such
that $\tau=\sigma\cap\Ker\,u$; the linear form $u$ descends to $\bar
u\in\bar\sigma^\vee$, therefore $\bar\sigma\cap\Ker\,\bar u=\{0\}$
is a face of $\bar\sigma$.

(ii): If $\tau$ has codimension two, $\bar V$ has dimension two.
Suppose that the assertion is known for $(\bar V,\bar\sigma)$; then
we find exactly two facets $\bar\gamma_1,\bar\gamma_2$ of
$\bar\sigma$ whose intersection is $\{0\}$. Their preimages in $V$
intersect $\sigma$ in facets $\gamma_1$, $\gamma_2$ that satisfy
$\gamma_1\cap\gamma_2=\tau$. Hence, we may assume from start that
$\tau=\{0\}$ and $V$ has dimension two, in which case the
verification is easy, and shall be left to the reader.

(i): Arguing by induction on the codimension, it suffices to show
that $\tau$ is contained in a proper face spanning a larger
subspace. To this aim, suppose that the claim is known for
$\bar\sigma$; since $\{0\}$ is a face of $\bar\sigma$ of codimension
at least two, it is contained in a proper face $\bar\gamma$; the
preimage $\gamma$ of the latter intersects $\sigma$ in a proper face
containing $\tau$. Thus again, we are reduced to the case where
$\tau=\{0\}$. Pick $u_0\in\sigma^\vee$ such that
$\sigma\cap\Ker\,u_0=\{0\}$; choose also any other $u_1\in V^\vee$
such that $\sigma\cap\Ker\,u_1\neq\{0\}$. Since $\dim_\R V^\vee\geq
2$, we may find a continuous map $f:[0,1]\to
V^\vee\!\setminus\!\{0\}$ with $f(0)=u_0$ and $f(1)=u_1$. Let
$\P_+(V)$ be the topological space of rays of $V$ ({\em i.e.} the
topological quotient $V\!\setminus\!\{0\}/\!\!\sim$ by the
equivalence relation such that $v\sim v'$ if and only if $v$ and
$v'$ generate the same ray), and define likewise $\P_+(V^\vee)$; let
also $Z\subset P:=\P_+(V)\times\P_+(V^\vee)$ be the incidence
correspondence, {\em i.e.} the subset of all pairs $(\bar v,\bar u)$
such that $u(v)=0$, for any representative $u$ of the class $\bar u$
and $v$ of the class $\bar v$. Finally, let
$\P_+(\sigma)\subset\P_+(V)$ be the image of
$\sigma\!\setminus\!\{0\}$. Then $Z$ (resp. $\P_+(\sigma)$) is a
closed subset of $P$ (resp. of $\P_+(V)$), hence
$Y:=Z\cap(\P_+(\sigma)\times\P_+(V^\vee))$ is a closed subset of
$P$. Since the projection $\pi:P\to\P_+(V^\vee)$ is proper, $\pi(Y)$
is closed in $\P_+(V^\vee)$. Let $\bar f:[0,1]\to\P_+(V^\vee)$ be
the composition of $f$ and the natural projection
$V^\vee\!\setminus\!\{0\}\to\P_+(V^\vee)$; then $f^{-1}(\pi(Y))$ is
a closed subset of $[0,1]$, hence it admits a smallest element, say
$a$ (notice that $a>0$). Moreover, $u_a\in\sigma^\vee$; indeed,
otherwise we may find $v\in\sigma\!\setminus\!\{0\}$ such that
$u_a(v)<0$, and since $u_0(v)>0$, we would have $u_b(v)=0$ for some
$b\in(0,a)$. The latter means that $f(b)\in\pi(Y)$, which
contradicts the definition of $a$. Since by construction,
$\sigma\cap\Ker\,u_a\neq\{0\}$, the claim follows.
\end{pfclaim}

Let $\tau$ be any face of $\sigma$; to show that (iii) holds for
$\tau$, we argue by induction on the codimension of $\tau$. The case
of codimension 2 is covered by claim \ref{cl_codim-face}(ii). If
$\tau$ has codimension $>2$, we apply claim \ref{cl_codim-face}(i)
to find a proper face $\gamma$ containing $\tau$; by induction,
$\tau$ is the intersection of facets of $\gamma$, and each of these
is the intersection of two facets in $\sigma$ (again by claim
\ref{cl_codim-face}(ii)), whence the contention.
\end{proof}

\sset\subsubsection{}\label{subsec_spann}
Suppose $\sigma$ spans $V$ ({\em i.e.} $\langle\sigma\rangle=V$)
and let $\tau$ be a facet of $\sigma$; by definition there exists
an element $u_\tau\in\sigma^\vee$ such that
$\tau=\sigma\cap\Ker\,u_\tau$, and one sees easily that the ray
$R_\tau:=\R_+\cdot u_\tau\subset\sigma^\vee$ is well-defined,
independently of the choice of $u_\tau$. Hence the half-space :
$$
H_\tau:=\{v\in V~|~u_\tau(v)\geq 0\}
$$
depends only on $\tau$. Recall that the {\em interior\/} (resp.
the {\em topological closure}) of a subset $E\subset V$ is the
largest open subset (resp. the smallest closed subset) of $V$
contained in $E$ (resp. containing $E$). The {\em topological
boundary\/} of $E$ is the intersection of the topological
closures of $E$ and of its complement $V\!\setminus\!E$.

\begin{proposition}\label{prop_spanning-cone}
Let $(V,\sigma)$ be a convex polyhedral cone, such that $\sigma$
spans $V$. We have:
\begin{enumerate}
\item
The topological boundary of $\sigma$ is the union of its facets.
\item
If $\sigma\neq V$, then $\sigma=\bigcap_{\tau\subset\sigma}H_\tau$,
where $\tau$ ranges over the facets of $\sigma$.
\item
The rays $R_\tau$, where $\tau$ ranges over the facets of $\sigma$,
generate the cone $\sigma^\vee$.
\end{enumerate}
\end{proposition}
\begin{proof} (i): Notice that $\sigma$ spans $V$ if and only if
the interior of $\sigma$ is not empty. A proper face $\tau$ is the
intersection of $\sigma$ with a hyperplane $\Ker\,u\subset V$ with
$u\in\sigma^\vee\!\setminus\!\{0\}$; therefore, every neighborhood
$U\subset V$ of any point $v\in\tau$ intersects
$V\!\setminus\!\sigma$. This shows that $\tau$ lies in the
topological boundary of $\sigma$.

Conversely, if $v$ is in the boundary of $\sigma$, choose a sequence
$(v_i~|~i\in\N)$ of points of $V\!\setminus\!\sigma$, converging to
the point $v$; by lemma \ref{lem_ddoble-cone}, for every $i\in\N$
there exists $u_i\in\sigma^\vee$ such that $u_i(v_i)<0$. Up to
multiplication by scalars, we may assume that the vectors $u_i$ lie
on some sphere in $V^\vee$ (choose any norm on $V^\vee$); hence we
may find a convergent subsequence, and we may then assume that the
sequence $(u_i~|~i\in\N)$ converges to an element $u\in V^\vee$.
Necessarily $u\in\sigma^\vee$, therefore $v$ lies on the face
$\sigma\cap\Ker\,u$, and the assertion follows from proposition
\ref{prop_refer-to-Ful}(iii).

(ii): Suppose, by way of contradiction, that $v$ lies in every
half-space $H_\tau$, but $v\notin\sigma$. Choose any point $v'$ in
the interior of $\sigma$, and let $t\in[0,1]$ be the largest value
such that $w:=tv+(1-t)v'\in\sigma$. Clearly $w$ lies on the boundary
of $\sigma$, hence on some facet $\tau$, by (i). Say that
$\tau=\sigma\cap\Ker\,u$; then $u(v')>0$ and $u(w)=0$, so $u(v)<0$,
a contradiction.

(iii): When $\sigma=V$, there is nothing to prove, hence we may
assume that $\sigma\neq V$. In this case, suppose that
$u\in\sigma^\vee$, and $u$ is not in the cone $C$ generated by the
rays $R_\tau$. Applying lemma \ref{lem_ddoble-cone} to the cone
$(V^\vee,C)$, we deduce that there exists a vector $v\in V$ with
$v\in H_\tau$ for all the facets $\tau$ of $\sigma$, and $u(v)<0$,
which contradicts (ii).
\end{proof}

\begin{corollary}\label{cor_spanning-cone}
Let $(V,\sigma)$ and $(V,\sigma')$ be two convex polyhedral cones.
Then :
\begin{enumerate}
\item
{\em (Farkas' theorem)}\ \ The dual $(V^\vee,\sigma^\vee)$ is also a
convex polyhedral cone.
\item
If $\tau$ is a face of\/ $\sigma$, then
$\tau^*:=\sigma^\vee\cap\tau^\bot$ is a face of\/ $\sigma^\vee$ such
that $\langle\tau^*\rangle=\langle\tau\rangle^\bot$. Especially:
\set\begin{equation}\label{eq_dim-of-cplt}
\dim_\R\langle\tau\rangle+\dim_\R\langle\tau^*\rangle=\dim_\R V.
\end{equation}
The rule $\tau\mapsto\tau^*$ is a bijection from the set of faces
of\/ $\sigma$ to those of\/ $\sigma^\vee$. The smallest face of\/
$\sigma$ is $\sigma\cap(-\sigma)$.
\item
$(V,\sigma\cap\sigma')$ is a convex polyhedral cone, and every face
of\/ $\sigma\cap\sigma'$ is of the form $\tau\cap\tau'$, for some
faces $\tau$ of\/ $\sigma$ and $\tau'$ of\/ $\sigma'$.
\end{enumerate}
\end{corollary}
\begin{proof} (i): Set $W:=\langle\sigma\rangle\subset V$, and pick a
basis $u_1,\dots,u_k$ of $W^\bot$; by proposition
\ref{prop_spanning-cone}(iii), the assertion holds for the dual
$(W^\vee,\sigma^\vee)$ of the cone $(W,\sigma)$. However,
$W^\vee\simeq V^\vee/W^\bot$, hence the dual cone
$(V^\vee,\sigma^\vee)$ is generated by lifts of generators of
$(W^\vee,\sigma^\vee)$, together with the vectors $u_i$ and $-u_i$,
for $i=1,\dots,k$.

(ii): Notice first that the faces of $\sigma^\vee$ are exactly the
cones $\sigma^\vee\cap\{u\}^\bot$, for
$u\in\sigma=(\sigma^\vee)^\vee$. For a given $v\in\sigma$, let
$\tau$ be the smallest face of $\sigma$ such that $v\in\tau$; this
means that $\tau^\vee\cap\{v\}^\bot=\tau^\bot$ (where
$(V^\vee,\tau^\vee)$ is the dual of $(V,\tau)$). Hence
$\sigma^\vee\cap\{v\}^\bot=\sigma^\vee\cap(\tau^\vee\cap\{v\}^\bot)=\tau^*$,
so every face of $\sigma^\vee$ has the asserted form. The rule
$\tau\mapsto\tau^*$ is clearly order-reversing, and from the obvious
inclusion $\tau\subset(\tau^*)^*$ it follows that
$\tau^*=((\tau^*)^*)^*$, hence this map is a bijection. It follows
from this, that the smallest face of $\sigma$ is
$(\sigma^\vee)^*=\sigma\cap(\sigma^\vee)^\bot=(\sigma^\vee)^\bot=
\sigma\cap(-\sigma)$. In particular, we see that
$(\sigma\cap(-\sigma))^*=\sigma^\vee$, and furthermore,
\eqref{eq_dim-of-cplt} holds for $\tau:=\sigma\cap(-\sigma)$.
Identity \eqref{eq_dim-of-cplt} for a general face $\tau$ can be
deduced by inserting $\tau$ in a maximal chain of faces of $\sigma$,
and comparing with the dual chain of faces of $\sigma^\vee$ (details
left to the reader). Finally, it is clear that
$\langle\tau\rangle\subset\langle\tau^*\rangle^\bot$; since these
spaces have the same dimension, we deduce that
$\langle\tau\rangle^\bot=\langle\tau^*\rangle$.

(iii): Indeed, lemma \ref{lem_ddoble-cone} implies that
$\sigma^\vee+\sigma^{\prime\vee}$ is the dual of
$\sigma\cap\sigma'$, hence (i) implies that $\sigma\cap\sigma'$ is
polyhedral. It also follows that every face $\tau$ of
$\sigma\cap\sigma'$ is the intersection of $\sigma\cap\sigma'$ with
the kernel of a linear form $u+u'$, for some $u\in\sigma^\vee$ and
$u'\in\sigma^{\prime\vee}$. Consequently,
$\tau=(\sigma\cap\Ker\,u)\cap(\sigma'\cap\Ker\,u')$.
\end{proof}

\begin{corollary}\label{cor_strongly}
For a convex polyhedral cone $(V,\sigma)$, the following conditions
are equivalent:
\begin{enumerate}
\alphaenu
\item
$\sigma\cap(-\sigma)=\{0\}$.
\item
$\sigma$ contains no non-zero linear subspaces.
\item
There exists $u\in\sigma^\vee$ such that $\sigma\cap\Ker\,u=\{0\}$.
\item
$\sigma^\vee$ spans $V^\vee$.
\end{enumerate}
\end{corollary}
\begin{proof} (a) $\Leftrightarrow$ (b) since $\sigma\cap(-\sigma)$
is the largest subspace contained in $\sigma$. Next, (a)
$\Leftrightarrow$ (c) since $\sigma\cap(-\sigma)$ is the smallest
face of $\sigma$. Finally, (a) $\Leftrightarrow$ (d) since
$\dim_\R\langle\sigma\cap(-\sigma)\rangle+
\dim_\R\langle\sigma^\vee\rangle=n$ (corollary
\ref{cor_spanning-cone}(ii)).
\end{proof}

\sset\subsubsection{}\label{subsec_extremal}
A convex polyhedral cone fulfilling the equivalent conditions of
corollary \ref{cor_strongly} is said to be {\em strictly convex}.
Suppose that $(V,\sigma)$ is strictly convex; then proposition
\ref{prop_spanning-cone}(iii) says that $\sigma$ is generated by the
rays $R_\tau$, where $\tau$ ranges over the facets of $\sigma^\vee$.
The rays $R_\tau$ are uniquely determined by $\sigma$, and are
called the {\em extremal rays} of $\sigma$. Moreover, these $R_\tau$
form {\em the unique minimal set of generating rays\/} for $\sigma$.
Indeed, concerning the minimality : for each facet $\tau$ of
$\sigma^\vee$, pick $v_\tau\in\sigma$ with $\R_+\cdot
v_\tau=R_\tau$; suppose that $v_{\tau_0}=\sum_{\tau\in S}t_\tau\cdot
v_\tau$ for some subset $S$ of the set of facets of $\sigma^\vee$,
and appropriate $t_\tau>0$, for every $\tau\in S$. It follows easily
that $u(v_\tau)=0$ for every $u\in\tau_0$, and every $\tau\in S$.
But by definition of $R_\tau$, this implies that $S=\{\tau_0\}$,
which is the claim. Concerning uniqueness : suppose that $\Sigma$ is
another system of generating rays; especially, for any facet
$\tau\subset\sigma^\vee$, the ray $R_\tau$ is in the convex cone
generated by $\Sigma$; it follows easily that there exists
$\rho\in\Sigma$ such that $u(\rho)=0$ for every $u\in\tau$, in which
case $\rho=R_\tau$. This shows that $\Sigma$ must contain all the
extremal rays.

\begin{example}\label{ex_cone-dim-two}
(i)\ \
Suppose that $\dim_\R V=2$, and $(V,\sigma)$ is a strictly convex
polyhedral cone, and assume that $\sigma$ generates $V$. Then the
only face of codimension two of $\sigma$ is $\{0\}$, so it follows
easily from claim \ref{cl_codim-face}(ii) that $\sigma$ admits
exactly two facets, and these are also the extremal rays of $\sigma$,
especially $\sigma$ is simplicial. Of course, these assertions are
rather obvious; in dimension $>2$, the general situation is much
more complicated.

(ii)\ \
Let $(V,\sigma)$ be a convex polyhedral cone, and suppose that
$\sigma$ spans $V$. Let $\tau$ be a face of $\sigma$. Notice
that $(\sigma,+)^\gp=(V,+)$, and $(\tau,+)$ is a face of the monoid
$(\sigma,+)$, by lemma \ref{lem_same-faces}. Hence we may
view the localization $\tau^{-1}\sigma$ naturally as a submonoid
of $(V,+)$, and it is easily seen that $\tau^{-1}\sigma$ is
a convex cone. By proposition \ref{prop_spanning-cone}(iii),
the polyhedral cone $\sigma^\vee$ is generated by its extremal
rays $\R l_1,\dots,\R l_n$, and by proposition
\ref{prop_refer-to-Ful}(iii), we may assume that
$\tau=\sigma\cap\Ker\,(l_1+\cdots+l_k)$ for some $k\leq n$.
Clearly $l_i(v)\geq 0$ for every $v\in\tau^{-1}\sigma$ and
every $i\leq k$. Conversely, if $l\in(\tau^{-1}\sigma)^\vee$,
we must have $\tau\subset\Ker\,l$ and $l\in\sigma^\vee$;
if we write $l=\sum_{i=1}^na_il_i$ for some $a_i\geq 0$,
it follows easily that $a_i=0$ for every $i>k$. On the other
hand, suppose that $v\in V$ satisfies the inequalities
$l_i(v)\geq 0$ for $i=1,\dots,k$; then, for every
$i=k+1,\dots,n$ we may find $u_i\in\tau$ such that
$l_i(v+u_i)\geq 0$, hence $v+u_{k+1}+\cdots+u_n\in\sigma$,
and therefore $v\in\tau^{-1}\sigma$. This shows that
$\tau^{-1}\sigma$ is a closed convex cone, and its dual
$(\tau^{-1}\sigma)^\vee$ is the convex cone generated by
$l_1,\dots,l_k$; especially, it is a convex polyhedral cone,
and then the same holds for $\tau^{-1}\sigma$, by virtue of
lemma \ref{lem_ddoble-cone} and corollary \ref{cor_spanning-cone}(i).

(iii)\ \
In the situation of (ii), let $v\in\tau$ be any element that
lies in the {\em relative interior\/} of $\tau$, {\em i.e.}
in the complement of the union of the facets of $\tau$. Denote
by $S_v\subset\tau$ the submonoid generated by $v$. Then we
claim that
$$
\tau^{-1}\sigma=S^{-1}_v\sigma.
$$
Indeed, let $s\in\sigma$ and $t\in\tau$ be any two element;
in view of proposition \ref{prop_spanning-cone}(i),
it is easily seen that there exists an integer $N>0$ such
that $v-N^{-1}t\in\tau$, hence $t':=Nv-t\in\tau$. Therefore
$s-t=(s+t')-Nv\in S^{-1}\sigma$, and the assertion follows.
\end{example}

\begin{lemma}\label{lem_dir-img-cone}
Let $f:V\to W$ be a linear map of finite dimensional $\R$-vector
spaces, $(V,\sigma)$ a convex polyhedral cone. The following holds :
\begin{enumerate}
\item
$(W,f(\sigma))$ is a convex polyhedral cone.
\item
Suppose moreover, that $\sigma\cap\Ker\,f$ does not contain non-zero
linear subspaces of $V$. Then, for every face $\tau$ of $f(\sigma)$
there exists a face $\tau'$ of $\sigma$ such that
$f(\tau')\subset\tau$, and $f$ restricts to an isomorphism :
$\langle\tau'\rangle\isom\langle\tau\rangle$.
\end{enumerate}
\end{lemma}
\begin{proof}(i) is obvious. To show (ii) we argue by induction
on $n:=\dim_\R\Ker\,f$. The assertion is obvious when $n=0$, hence
suppose that $n>0$ and that the claim is already known whenever
$\Ker\,f$ has dimension $<n$. Let $\tau$ be a face of $f(\sigma)$;
then $f^{-1}\tau$ is a face of $f^{-1}f(\sigma)$, hence $\sigma\cap
f^{-1}\tau$ is a face of $\sigma=\sigma\cap f^{-1}f(\sigma)$. In
view of proposition \ref{prop_refer-to-Ful}(ii), we may then replace
$\sigma$ by $\sigma\cap f^{-1}\tau$, and therefore assume from start
that $\tau=f(\sigma)$. We may as well assume that
$V=\langle\sigma\rangle$ and $W=\langle\tau\rangle$. The assumption
on $\sigma$ implies especially that $\sigma\neq V$, hence $\sigma$
is the intersection of the half-spaces $H_\gamma$ corresponding to
its facets $\gamma$ (proposition \ref{prop_spanning-cone}(ii)). For
each facet $\gamma$ of $\sigma$, let $u_\gamma$ be a chosen
generator of the ray $R_\gamma$ (notation of \eqref{subsec_spann}).
Since $\sigma\cap\Ker\,f$ does not contain non-zero linear
subspaces, we may find a facet $\gamma$ such that
$V':=\Ker\,u_\gamma$ does not contain $\Ker\,f$. Then, the inductive
assumption applies to the restriction $f_{|V'}:V'\to W$ and the
convex polyhedral cone $\sigma\cap V'$, and yields a face $\tau'$ of
the latter, such that $f$ induces an isomorphism
$\langle\tau'\rangle\isom\langle f(\sigma\cap V')\rangle$. Finally,
$\langle f(\sigma\cap V')\rangle=W$, since $V'$ does not contain
$\Ker\,f$.
\end{proof}

\begin{lemma}\label{lem_pull-a-a-cone}
Let $f:V\to W$ be a linear map of finite dimensional $\R$-vector
spaces, $f^\vee:W^\vee\to V^\vee$ the transpose of $f$, and
$(W,\sigma)$ a convex polyhedral cone. Then :
\begin{enumerate}
\item
$(V,f^{-1}\sigma)$ is a convex polyhedral cone and
$(f^{-1}\sigma)^\vee=f^\vee(\sigma^\vee)$.
\item
For every face $\delta$ of $f^{-1}\sigma$, there exists a face
$\tau$ of $\sigma$ such that $\delta=f^{-1}\tau$. If furthermore,
$\langle\sigma\rangle+f(V)=W$, we may find such a $\tau$ so that
additionally, $f$ induces an isomorphism:
$$
V/\langle\delta\rangle\isom W/\langle\tau\rangle.
$$
\item
Conversely, for every face $\tau$ of $\sigma$, the cone $f^{-1}\tau$
is a face of $f^{-1}\sigma$, and $(f^{-1}\tau)^*$ is the smallest
face of $(f^{-1}\sigma)^\vee$ containing $f^\vee(\tau^*)$
(notation of corollary {\em\ref{cor_spanning-cone}(ii)}).
\end{enumerate}
\end{lemma}
\begin{proof} (i): By corollary \ref{cor_spanning-cone}(i),
$(W^\vee,\sigma^\vee)$ is a convex polyhedral cone, hence we may
find $u_1,\dots,u_s\in\sigma^\vee$ such that
$\sigma=\bigcap_{i=1}^su_i^{-1}(\R_+)$. Therefore
$f^{-1}\sigma=\bigcap_{i=1}^s(u_i\circ f)^{-1}(\R_+)$. Let
$\gamma\subset V^\vee$ be the cone generated by the set
$\{u_1\circ f,\dots,u_s\circ f\}$; then $f^{-1}\sigma=\gamma^\vee$,
and the assertion results from lemma \ref{lem_ddoble-cone} and
a second application of (i).

(iii): For every $u\in\sigma^\vee$, we have :
$f^{-1}(\sigma\cap\Ker\,u)=(f^{-1}\sigma)\cap\Ker\,u\circ f$. Since
we already know that $(f^{-1}\sigma)^\vee=f^\vee(\sigma^\vee)$, we
see that the faces of $f^{-1}\sigma$ are exactly the subsets of the
form $f^{-1}\tau$, where $\tau$ ranges over the faces of $\sigma$.
Next, for any such $\tau$, the set $f^\vee(\tau^*)$ consists of all
$u\in V^\vee$ of the form $u=w\circ f$ for some $w\in\sigma^\vee$
such that $w(\tau)=0$. From this description it is clear that
$f^\vee(\tau^*)\subset(f^{-1}\tau)^*$. To show that $(f^{-1}\tau)^*$
is the smallest face containing $f^\vee(\tau^*)$, it then suffices
to prove that $f^\vee(\tau^*)^\bot\cap f^{-1}\sigma\subset
f^{-1}\tau$. However, let $v\in f^\vee(\tau^*)^\bot\cap
f^{-1}\sigma$; then $w\circ f(v)=0$ for every $w\in\tau^*$, {\em
i.e.} $f(v)\in(\tau^*)^*=\tau$, whence the contention.

(ii): The first assertion has already been shown; hence, suppose
that $\langle\sigma\rangle+f(V)=W$. We deduce that
$\sigma^\vee\cap\Ker\,f^\vee$ does not contain non-zero linear
subspaces of $W^\vee$; indeed, if $u\in W^\vee$, and both $u$ and
$-u$ lie in $\sigma^\vee$, then $u$ vanishes on
$\langle\sigma\rangle$, and if $u\in\Ker\,f^\vee$, then $u$ vanishes
as well on $f(V)$, hence $u=0$. We may then apply lemma
\ref{lem_dir-img-cone}(ii) to find a face $\gamma$ of $\sigma^\vee$
such that $f^\vee(\gamma)\subset\delta^*$ and $f^\vee$ restricts to
an isomorphism : $\langle\gamma\rangle\isom\langle\delta^*\rangle$.
Especially, $\delta^*$ is the smallest face of $(f^{-1}\sigma)^\vee$
containing $f^\vee(\gamma)$, hence $\delta^*=(f^{-1}\gamma^*)^*$ by
(iii), {\em i.e.} $\delta=f^{-1}\gamma^*$. We also deduce that $f$
induces an isomorphism : $V/\langle\delta^*\rangle^\bot\isom
W/\langle\gamma\rangle^\bot$. Since
$\langle\delta^*\rangle^\bot=\langle\delta\rangle$ and
$\langle\gamma\rangle^\bot=\langle\gamma^*\rangle$, the second
assertion holds with $\tau:=\gamma^*$.
\end{proof}

\begin{lemma}\label{lem_prod-of-polyh}
Let $(V,\sigma)$ and $(V',\sigma')$ be two convex polyhedral
cones. Then :
\begin{enumerate}
\item
$(V\oplus V',\sigma\times\sigma')$ is a convex polyhedral cone.
\item
Every face of $\sigma\times\sigma'$ is of the form
$\tau\times\tau'$, for some faces $\tau$ of $\sigma$ and $\tau'$ of
$\sigma'$.
\end{enumerate}
\end{lemma}
\begin{proof} Indeed,
$\sigma\times\sigma'=(p^{-1}_1\sigma)\cap(p^{-1}_2\sigma')$, where
$p_1$ and $p_2$ are the natural projections of $V\oplus V'$ onto $V$
and $V'$. Hence, assertions (i) and (ii) follow from corollary
\ref{cor_spanning-cone}(iii) and lemma \ref{lem_pull-a-a-cone}(i).
\end{proof}

\sset\subsubsection{}\label{subsec_from-con-to-mon}
Let $(L,+)$ be a free abelian group of finite rank, $\sigma\subset
L_\R:=L\otimes_\Z\R$ a convex polyhedral cone. We say that $\sigma$
is {\em $L$-rational\/} (or briefly : {\em rational\/}, when there
is no danger of ambiguity) if it admits a generating set consisting
of elements of $L$. Then it is clear that every face of a rational
convex polyhedral cone is again rational (see
\eqref{subsec_rather-see-here}).  On the other hand, let
$$
(M,+)\subset(L,+)
$$
be a submonoid of $L$; we shall denote by $(L_\R,M_\R)$ the convex
cone generated by $M$ ({\em i.e.} the smallest convex cone in $L_\R$
containing the image of $M$). If $M$ is fine, $M_\R$ is a convex
polyhedral cone. Later we shall also find useful to consider the
subset :
$$
M_\Q:=\{m\otimes q~|~m\in M, q\in\Q_+\}\subset L_\Q:=L\otimes_\Z\Q
$$
which is a submonoid of $L_\Q$.

\begin{proposition}\label{prop_was-part-of-Gordon}
Let $(L,+)$ be a free abelian group of finite rank, with dual
$$
L^\vee:=\Hom_\Z(L,\Z).
$$
Let also $(L_\R,\sigma)$ and $( L_\R,\sigma')$ be two $L$-rational
convex polyhedral cones. We have :
\begin{enumerate}
\item
The dual $(L^\vee_\R,\sigma^\vee)$ is an $L^\vee$-rational convex
polyhedral cone.
\item
$( L_\R,\sigma\cap\sigma')$ is also an $L$-rational convex
polyhedral cone.
\item
Let $g:L'\to L$ (resp. $h:L\to L'$) be a map of free abelian
groups of finite rank, and denote by $g_\R:L'_\R\to L_\R$
(resp. $h_\R:L_\R\to L'_\R$) the induced $\R$-linear map.
Then, $(L'_\R,g_\R^{-1}\sigma)$ and $(L'_\R,h_\R(\sigma))$
are $L'$-rational.
\item
Let $L'$ be another free abelian group of finite rank, and
$(L'_\R,\sigma')$ an $L'$-rational convex polyhedral cone.
Then $(L_\R\oplus L'_\R,\sigma\times\sigma')$ is $L\oplus L'$-rational.
\end{enumerate}
\end{proposition}
\begin{proof}(i) and (ii) follow easily, by inspecting the proof of
corollary \ref{cor_spanning-cone}(i),(iii) : the details shall be
left to the reader.

(iii): The assertion concerning $h_\R(\sigma)$ is obvious.
To show the assertion for $g^{-1}_\R\sigma$, one argues as
in the proof of lemma \ref{lem_pull-a-a-cone}(i) : by (i),
we may find $u_1,\dots,u_s\in L^\vee$ such that
$\sigma=\bigcap_{i=1}^su_{i,\R}^{-1}(\R_+)$. Therefore
$g_\R^{-1}\sigma=\bigcap_{i=1}^s(u_i\circ g)_\R^{-1}(\R_+)$. Let
$\gamma\subset V^\vee$ be the cone generated by the set
$\{(u_1\circ g)_\R,\dots,(u_s\circ g)_\R\}$; then
$g_\R^{-1}\sigma=\gamma^\vee$, and the assertion results from
lemma \ref{lem_ddoble-cone} and a second application of (i).

Lastly, arguing as in the proof of lemma \ref{lem_prod-of-polyh}(i),
one derives (iii) from (ii) and (iii).
\end{proof}

Parts (i) and (iii) of the following proposition provide the
bridge connecting convex polyhedral cones to fine monoids.

\begin{proposition}\label{prop_Gordon}
Let $(L,+)$ be a free abelian group of finite rank,
$(L_\R,\sigma)$ an $L$-rational convex polyhedral cone, and
set $\sigma_L:=L\cap\sigma$. Then :
\begin{enumerate}
\item
{\em (Gordan's lemma)}\ \ $\sigma_L$ is a fine and saturated
submonoid of $L$, and $L\cap\La\sigma\Ra=\sigma_L^\gp$.
\item
For every $v\in L_\R$, the subset $L\cap(\sigma-v)$ is a finitely
generated $\sigma_L$-module.
\item
For any submonoid $M\subset L$, we have : $M_\Q=M_\R\cap L_\Q$ and
$M^\sat=M_\R\cap L$.
\end{enumerate}
\end{proposition}
\begin{proof} (Here $\sigma-v\subset L_\R$ denotes the
translate of $\sigma$ by the vector $-v$, {\em i.e.}
the subset of all $w\in L_\R$ such that $w+v\in\sigma$.)
Choose $v_1,\dots,v_s\in L$ that generate $\sigma$, and set
$$
C_\eps:=\biggl\{
\sum_{i=1}^st_iv_i~|~\text{$t_i\in[0,\eps]$\ \ for $i=1,\dots,s$}
\biggr\}
\qquad
\text{for every $\eps>0$}.
$$

(i): Clearly $L\cap\sigma$ is saturated. Since $C_1$ is compact
and $L$ is discrete, $C_1\cap L$ is a finite set. We claim that
$C_1\cap L$ generates the monoid $\sigma_L$. Indeed, if
$v\in\sigma_L$, write $v=\sum_{i=1}^sr_iv_i$, with $r_i\geq 0$
for every $i=1,\dots,s$; hence $r_i=m_i+t_i$ for some $m_i\in\N$
and $t_i\in[0,1[$, and therefore $v=v'+\sum_{i=1}^sm_iv_i$, where
$v',v_1,\dots,v_s\in C_1\cap L$. Next, it is clear that
$\sigma_L^\gp\subset L\cap\La\sigma\Ra$; for the converse,
say that $w\in L\cap\La\sigma\Ra$, and write $w=w_1-w_2$,
for some $w_1,w_2\in\sigma$. Then $w_1=\sum_{i=1}^st_iv_i$
for some $t_i\geq 0$; we pick $t'_i\in\N$ such that
$t'_i\geq t_i$ for every $i\leq s$, and we set
$w'_i:=\sum_{i=1}^st'_iv_i$. It follows that $w=w'_1-w'_2$,
where $w'_2:=w_2+(w'_1-w_1))$, and notice that $w'_1\in\sigma_L$
and $w'_2\in\sigma$; then we must have $w'_2\in\sigma_L$ as
well, and therefore $w\in\sigma^\gp_L$.

(ii) is similar : from the compactness of $C_1$ one sees that
$L\cap(C_1-v)$ is a finite set; on the other hand, arguing
as in the proof of (i), one checks easily that the latter
set generates the $\sigma_L$-module $L\cap(\sigma-v)$.

(iii): Let $x\in M_\R\cap L_\Q$; then we may write
\set\begin{equation}\label{eq_like-this-x}
x=\sum^n_{i=1}r_i\otimes m_i
\qquad
\text{where $m_i\in M$, $r_i>0$ for every $i\leq n$}.
\end{equation}

\begin{claim}\label{cl_like-this-x}
In the situation of proposition \ref{prop_Gordon}(iii), let
$x\in M_\R\cap L_\Q$, and write $x$ as in \eqref{eq_like-this-x}.
Then, for every $\eps>0$ there exist $q_1,\dots,q_n\in\Q_+$ with
$|r_i-q_i|<\eps$ for every $i=1,\dots,n$, and such that
$x=\sum^n_{i=1}q_i\otimes m_i$.
\end{claim}
\begin{pfclaim} Up to a reordering, we may assume that
$m_1,\dots,m_k$ form a basis of the $\Q$-vector space
generated by $m_1,\dots,m_n$, therefore
$m_{k+i}=\sum_{j=1}^kq_{ij}m_j$ for a matrix
$$
A:=(q_{ij}~|~i=1,\dots,n-k;\ j=1,\dots,k)
$$
with entries in $\Q$. Let $\underline r:=(r_1,\dots,r_k)$ and
$\underline r':=(r_{k+1},\dots,r_n)$; since $x\in L_\Q$, we deduce
that $\underline b:=\underline r+\underline r'\cdot A\in\Q^{\oplus
k}$. Moreover, if $\underline s:=(s_1,\dots,s_k)\in\Q^{\oplus k}$
and $\underline s':=(s_{k+1},\dots,s_n)\in\Q^{\oplus n-k}$ satisfy
the identity $\underline b=\underline s+\underline s'\cdot A$, then
$\sum_{i=1}^ns_i\otimes m_i=x$. If we choose $\underline s'$ very
close to $\underline r'$, then $\underline s$ shall be very close to
$\underline r$; especially, we can achieve that both $\underline s$
and $\underline s'$ are vectors with positive coordinates.
\end{pfclaim}

Claim \ref{cl_like-this-x} shows that $x\in M_\Q$, whence the first
stated identity; for the second identity, we are reduced to showing
that $M^\sat=M_\Q\cap L$, which is immediate.
\end{proof}

For various algebraic and geometric applications of the theory
of polyhedral cones, one is led to study subdivisions of a given
cone, in the sense of the following definition \ref{def_fans}.
Later we shall see a more abstract notion of subdivision,
in the context of general fans, which however finds its roots
and motivation in the intuitive manipulations of polyhedra that
we formalize hereafter.

\begin{definition}\label{def_fans}
Let $V$ be a finite dimensional $\R$-vector space.
\begin{enumerate}
\item
A {\em fan\/} in $V$ is a finite set $\Delta$ consisting of
convex polyhedral cones of $V$, such that :
\begin{itemize}
\item
for every $\sigma\in\Delta$, and every face $\tau$ of $\sigma$, also
$\tau\in\Delta$;
\item
for every $\sigma,\tau\in\Delta$, the intersection $\sigma\cap\tau$
is also an element of $\Delta$, and is a face of both $\sigma$ and
$\tau$.
\end{itemize}
\item
We say that $\Delta$ is a {\em simplicial fan\/} if all the elements
of $\Delta$ are simplicial cones.
\item
Suppose that $V=L\otimes_\Z\R$ for some free abelian group $L$; then
we say that $\Delta$ is {\em $L$-rational\/} if the same holds for
every $\tau\in\Delta$.
\item
A {\em refinement\/} of the fan $\Delta$ is a fan $\Delta'$ in $V$
with $\bigcup_{\sigma\in\Delta}\sigma=\bigcup_{\tau\in\Delta'}\tau$,
and such that every $\tau\in\Delta$ is the union of the
$\gamma\in\Delta'$ contained in $\tau$.
\item
A {\em
subdivision\/} of a convex polyhedral cone $(V,\sigma)$ is a
refinement of the fan $\Delta_\sigma$ consisting of $\sigma$
and its faces.
\end{enumerate}
\end{definition}

\begin{lemma}\label{lem_lem_top-dim-subdi}
Let $(V,\sigma)$ be any convex polyhedral cone, $\Delta$
a subdivision of $(V,\sigma)$. We let
$$
\Delta^s:=\{\tau\in\Delta~|~\La\tau\Ra=\La\sigma\Ra\}.
$$
Then $\bigcup_{\tau\in\Delta^s}\tau=\sigma$.
\end{lemma}
\begin{proof} Let $\tau_0\in\Delta$ be any element.
Then $\sigma':=\bigcup_{\tau\neq\tau_0}\tau$ is a closed
subset of $\sigma$. If $\La\tau_0\Ra\neq\La\sigma\Ra$, then
$\sigma\setminus\tau_0$ is a dense open subset of $\sigma$
contained in $\sigma'$; it follows that $\sigma'=\sigma$
in this case. Especially,
$\tau_0=\bigcup_{\tau\neq\tau_0}(\tau\cap\tau_0)$; since
each $\tau\cap\tau_0$ is a face of $\tau_0$, we see that
there must exist $\tau\neq\tau_0$ such that $\tau_0$ is
a face of $\tau$. The lemma follows immediately.
\end{proof}

\begin{example}\label{ex_roof-again}
(i)\ \
Certain useful subdivisions of a polyhedral cone $\sigma$
are produced by means of auxiliary real-valued functions
defined on $\sigma$. Namely, let us say that a continuous
function $f:\sigma\to\R$ is a {\em roof}, if the following
holds. There exist finitely many $\R$-linear forms
$l_1,\dots,l_n$ on $V$, such that
$f(v)=\min(l_1(v),\dots,l_n(v))$ for every $v\in\sigma$.
The concept of roof shall be reintroduced in section
\ref{sec_special-sub}, in a more abstract and general guise;
however, in order to grasp the latter, it is useful to keep
in mind its more concrete polyhedral incarnation.
We attach to $f$ a subdivision of $\sigma$, as follows.
For every $i,j=1,\dots,n$ define
$l_{ij}:=l_i-l_j$, and let $\tau_i\subset V^\vee$ be the
polyhedral cone $\sigma^\vee+\R l_{i1}+\cdots+\R l_{in}$.
From the identity $l_{ik}=l_{ij}+l_{jk}$ we easily deduce
that $\tau_i^\vee\cap\tau_j^\vee$ is a face of both
$\tau^\vee_i$ and $\tau_j^\vee$, for every $i,j=1,\dots,n$.
Denote by $\Theta$ the smallest subdivision of $\sigma$
containing all the $\tau^\vee_i$; it is easily seen that
$$
\sigma=\bigcup_{i=1}^n\tau^\vee_i
$$
and the restriction of $f$ to each $\tau_i^\vee$ agrees
with $l_i$.

(ii)\ \
Conversely, let $f:\sigma\to\R$ a continuous function; suppose
there exist a subdivision $\Theta$ of $\sigma$, and a system
$(l_\tau~|~\tau\in\Theta)$ of $\R$-linear forms on $V$ such
that
\begin{itemize}
\item
$f(v)=l_\tau(v)$ for every $\tau\in\Theta$ and every $v\in\tau$.
\item
$f(u+v)\geq f(u)+f(v)$ for every $u,v\in\sigma$.
\end{itemize}
Then we claim that $f$ is a roof on $\sigma$. Indeed, let
$\Theta^s\subset\Theta$ be the subset of all $\tau$ that span
$\La\sigma\Ra$. Notice fist that the system
$(l_\tau~|~\tau\in\Theta^s)$ already determines $f$ uniquely,
by virtue of lemma \ref{lem_lem_top-dim-subdi}. Next, let
$\tau,\tau'\in\Theta^s$ be any two elements, and pick an
element $v$ of the interior of $\tau$. For any $u\in\tau'$
and any $\eps>0$ we have, by assumption :
$f(v+\eps u)\geq f(v)+f(\eps u)$. If $\eps$ is small enough,
we have as well $v+\eps u\in\tau$, in which case the foregoing
inequality can be written as :
$$
l_\tau(v)+\eps\cdot l_\tau(u)=l_\tau(v+\eps u)\geq
l_\tau(v)+\eps\cdot l_{\tau'}(u)
$$
whence $l_\tau(u)\geq l_{\tau'}(u)=f(u)$ and the assertion
follows.
\end{example}

\begin{proposition}\label{prop_divide-et-imp}
Let $f:V\to W$ be a linear map of finite dimensional $\R$-vector
spaces, $(V,\sigma)$ a convex polyhedral cone, and
$h:\sigma\to f(\sigma)$ the restriction of $f$. Then :
\begin{enumerate}
\item
There exists a subdivision $\Delta$ of $(W,f(\sigma))$ such that :
$$
h^{-1}(a+b)=h^{-1}(a)+h^{-1}(b)
\qquad\text{for every $\tau\in\Delta$ and every $a,b\in\tau$}.
$$
\item
Suppose moreover that $V=L\otimes_\Z\R$, $W=L'\otimes_\Z\R$ and
$f=g\otimes_\Z\one_\R$ for a map $g:L\to L'$ of free abelian groups.
If $\sigma$ is $L$-rational, then we may find an $L$-rational
subdivision $\Delta$ such that {\em (i)} holds.
\end{enumerate}
\end{proposition}
\begin{proof} Let $V_0$ be the largest linear subspace contained
in $\sigma\cap\Ker\,f$. Notice that, under the assumptions of (ii),
we have : $V_0=\R\otimes_\Z\Ker\,g$. One verifies easily that the
proposition holds for the given map $f$ and the cone $(V,\sigma)$,
if and only if it holds for the induced map $\bar f:V/V_0\to
W/f(V_0)$ and the cone $(V/V_0,\bar\sigma)$ (where $\bar\sigma$ is
the image of $\sigma$ in $V/V_0$). Hence, we may replace $f$ by
$\bar f$, and assume from start that $\sigma\cap\Ker\,f$ contains no
non-zero linear subspaces. Moreover, we may assume that $\sigma$
spans $V$ and $f(\sigma)$ spans $W$.

(i): Let $S$ be the set of faces $\tau$ of $\sigma$ such that $f$
restricts to an isomorphism $\langle\tau\rangle\isom W$.

\begin{claim}\label{cl_extr-rays}
Let $\lambda\subset f(\sigma)$ be any ray. Then :
\begin{enumerate}
\item
$\lambda':=h^{-1}\lambda$ is a strictly convex polyhedral
cone. Especially, $\lambda'$ is generated by its extremal
rays (see \eqref{subsec_extremal}).
\item
For every extremal ray $\rho$ of $\lambda'$ with
$\rho\not\subset\Ker\,f$, there exists $\tau\in S$ such that
$\rho=\tau\cap f^{-1}(\lambda)$.
\end{enumerate}
\end{claim}
\begin{pfclaim}
$\lambda'$ is a convex polyhedral cone by lemma
\ref{lem_pull-a-a-cone}(i) and corollary
\ref{cor_spanning-cone}(iii). To see that $\lambda'$ is strictly
convex, notice that any subspace $L\subset f^{-1}(\lambda)$ lies
already in $\Ker\,f$, and if $L\subset\sigma$, we must have
$L=\{0\}$ by assumption. Let $\rho$ be an extremal ray of $\lambda'$
which is not contained in $\Ker\,f$; notice that $\lambda'$ is the
intersection of the polyhedral cones
$\lambda_1:=h^{-1}\langle\lambda\rangle$ and
$\lambda_2:=f^{-1}\lambda$, hence we can find faces $\delta_i$
of $\lambda_i$ ($i=1,2$) such that $\rho=\delta_1\cap\delta_2$
(corollary \ref{cor_spanning-cone}(iii) and lemma
\ref{lem_pull-a-a-cone}(i)). However, the only proper face
of $\lambda_2$ is $\Ker\,f$ (lemma \ref{lem_pull-a-a-cone}(ii)),
hence $\delta_2=\lambda_2$. Likewise, $f^{-1}\langle\lambda\rangle$
has no proper faces, hence $\delta_1=\gamma\cap
f^{-1}\langle\lambda\rangle$ for some face $\gamma$ of $\sigma$
(again by corollary \ref{cor_spanning-cone}(iii)). Since $\lambda_2$
is a half-space in $f^{-1}\langle\lambda\rangle$, we deduce easily
that either $\delta_1=\rho$ or $\delta_1=\langle\rho\rangle$.
Especially,
$\dim_\R(f^{-1}\langle\lambda\rangle)/\langle\delta_1\rangle=
\dim_\R\Ker\,f$. We may then apply lemma \ref{lem_pull-a-a-cone}(ii)
to the embedding $f^{-1}\langle\lambda\rangle\subset V$, to find a
face $\tau$ of $\sigma$ such that :
$$
\tau\cap f^{-1}\langle\lambda\rangle=\delta_1\qquad
\langle\tau\rangle\cap
f^{-1}\langle\lambda\rangle=\langle\rho\rangle \qquad\dim_\R
V/\langle\tau\rangle=\dim_\R\Ker\,f.
$$
It follows that $\langle\tau\rangle\cap\Ker\,f=\{0\}$, and therefore
$\tau\in S$, as required.
\end{pfclaim}

We construct as follows a subdivision of $(W,f(\sigma))$. For every
$\tau\in S$, let $F(\tau)$ be the set consisting of the facets of
the polyhedral cone $f(\tau)$; set also $F:=\bigcup_{\tau\in S}
F(\tau)$. Notice that, for every $\gamma\in F$, the subspace
$\langle\gamma\rangle$ is a hyperplane of $W$; we let :
$$
U:=f(\sigma)\!\setminus\!\bigcup_{\gamma\in F}\langle\gamma\rangle.
$$
Then $U$ is an open subset of $f(\sigma)$, and the topological
closure $\bar C$ of every connected component $C$ of $U$ is a convex
polyhedral cone. Moreover, if $C$ and $D$ are any two such connected
components, the intersection $\bar C\cap\bar D$ is a face of both
$\bar C$ and $\bar D$. We let $\Delta$ be the subdivision of
$f(\sigma)$ consisting of the cones $\bar C$ -- where $C$ ranges
over all the connected components of $U$ -- together with all their
faces.

\begin{claim}\label{cl_good-face}
For every $\delta\in\Delta$ and every $\tau\in S$, the intersection
$\delta\cap f(\tau)$ is a face of $\delta$.
\end{claim}
\begin{pfclaim} Due to proposition \ref{prop_refer-to-Ful}(iii), we may
assume that $\delta$ is the topological closure of a connected
component $C$ of $U$. We may also assume that $f(\tau)\neq W$,
otherwise there is nothing to prove; in that case, we have
$f(\tau)=\bigcap_{\gamma\in F(\tau)}H_\gamma$, where, for each
$\gamma\in F(\tau)$, the half-space $H_\gamma$ is the unique one
that contains both $f(\tau)$ and $\gamma$ (proposition
\ref{prop_spanning-cone}(ii)). It then suffices to show that
$\delta\cap H_\gamma$ is a face of $\delta$ for each such
$H_\gamma$. We may assume that $\delta\not\subset H_\gamma$. Since
$C$ is connected and $C\subset W\!\setminus\!\langle\gamma\rangle$,
it follows that $\delta\subset -H_\gamma$, the topological closure
of the complement of $H_\gamma$. Hence
$(-H_\gamma)^\vee\subset\delta^\vee$ (where $(-H_\gamma)^\vee$ is
the dual of the polyhedral cone $(W,-H_\gamma)$), and therefore
$\delta\cap H_\gamma=\delta\cap H_\gamma\cap(-H_\gamma)=
\delta\cap\langle\gamma\rangle$ is indeed a face of $\delta$.
\end{pfclaim}

Next, for every $w\in f(\sigma)$, let $I(w):=\{\tau\in S~|~w\in
f(\tau)\}$.

\begin{claim}\label{cl_intersect-I}
Let $\delta\in\Delta$, and $w_1,w_2\in\delta$. Then
$I(w_1+w_2)\subset I(w_1)\cap I(w_2)$.
\end{claim}
\begin{pfclaim} Suppose first that $w_1+w_2$ is contained in a
face $\delta'$ of $\delta$; say that $\delta'=\delta\cap\Ker\,u$,
for some $u\in\delta^\vee$. This means that $u(w_1+w_2)=0$, hence
$u(w_1)=u(w_2)=0$, {\em i.e.} $w_1,w_2\in\delta'$. Hence, we may
replace $\delta$ by $\delta'$, and assume that $\delta$ is the
smallest element of $\Delta$ containing $w_1+w_2$. Thus, suppose
that $\tau\in I(w_1+w_2)$; therefore $w_1+w_2\in f(\tau)\cap\delta$.
From claim \ref{cl_good-face} we deduce that $\delta\subset
f(\tau)$, hence $\tau\in I(w_1)\cap I(w_2)$, as claimed.
\end{pfclaim}

Finally, we are ready to prove assertion (i). Hence, let $a,b\in
f(\sigma)$ be any two vectors that lie in the same element of
$\Delta$. Clearly :
$$
h^{-1}(a)+h^{-1}(b)\subset h^{-1}(a+b)
$$
hence it suffices to show the converse inclusion. However, directly
from claim \ref{cl_extr-rays}(ii) we derive the identity :
$$
h^{-1}(\R_+\cdot w)= (\sigma\cap\Ker\,f)+
\sum_{\tau\in I(w)}(\tau\cap f^{-1}(\R_+\cdot w))
\qquad
\text{for every $w\in f(\sigma)$}.
$$
Taking into account claim \ref{cl_intersect-I}, we are then reduced
to showing that :
$$
\tau\cap f^{-1}(\R_+\cdot(a+b))\subset (\tau\cap f^{-1}(\R_+\cdot
a))+ (\tau\cap f^{-1}(\R_+\cdot b)) \qquad\text{for every $\tau\in
I(a+b)$}.
$$
The latter assertion is obvious, since $f$ restricts to an
isomorphism $\langle\tau\rangle\isom W$.

(ii): By inspecting the construction, one verifies easily that the
subdivision $\Delta$ thus exhibited shall be $L$-rational, whenever
$\sigma$ is.
\end{proof}

\sset\subsubsection{}\label{subsec_not-so-elegant}
Later we shall also be interested in rational variants of
the identities of proposition \ref{prop_divide-et-imp}(i).
Namely, consider the following situation.
Let $g:L\to L'$ be a map of free abelian groups of finite rank,
$g_\R:L_\R\to L'_\R$ the induced $\R$-linear map, and $(L_\R,\sigma)$
an $L$-rational convex polyhedral cone; set $\tau:=g_\R(\sigma)$,
and denote by $h_\R:\sigma\to\tau$ (resp.
$h_\Q:\sigma\cap L_\Q\to\tau\cap L'_\Q$) the restriction of $g_\R$.
We point out, for later reference, the following observation :

\begin{lemma}\label{lem_not-so-elegant}
In the situation of \eqref{subsec_not-so-elegant}, suppose that :
$$
h_\R^{-1}(x_1)+h_\R^{-1}(x_2)=h_\R^{-1}(x_1+x_2)
\qquad
\text{for every $x_1,x_2\in\tau$}
$$
(where the sum is taken in the monoid $(\cP(\sigma),+)$). Then
we have as well :
$$
h_\Q^{-1}(x_1)+h_\Q^{-1}(x_2)=h^{-1}_\Q(x_1+x_2)
\qquad
\text{for every $x_1,x_2\in\tau\cap L'_\Q$}.
$$
\end{lemma}
\begin{proof} Let $x_1,x_2\in\tau\cap L'_\Q$ be any two elements,
and $v\in h^{-1}_\Q(x_1+x_2)$, so we may write $v=v_1+v_2$ for
some $v_i\in h^{-1}_\R(x_i)$ ($i=1,2$). Let also $u_1,\dots,u_k$
be a finite system of generators for $\sigma^\vee$, and set
$$
J_i:=\{j\leq k~|~u_j(v_i)=0\}
\quad
E_i:=g_\R^{-1}(x_i)\cap\bigcap_{j\in J_i}\Ker\,u_j
\qquad
(i=1,2).
$$
Clearly $L_\Q\cap E_i$ is a dense subset of $E_i$ for $i=1,2$,
hence, in any neighborhood of $(x_1,x_2)$  in $L_\R^{\oplus 2}$
we may find a solution $(y_1,y_2)\in L^{\oplus 2}_\Q$ for the
system of equations
$$
g_\R(y_i)=x_i
\qquad
u_j(y_i)=0
\qquad
\text{for $i=1,2$ and every $j\in J_i$}.
$$
Since $u_j(x_i)>0$ for every $j\notin J_i$, we will also
have $u_j(y_i)>0$ for every $j\notin J_i$, provided $y_i$
is sufficiently close to $x_i$. The lemma follows.
\end{proof}

\sset\subsubsection{}\label{subsec_for-almost-lengths}
We conclude this section with some considerations that shall
be useful later, in our discussion of normalized lengths for
model algebras (see \eqref{subsec_tenero-giaco}). Keep the
notation of proposition \ref{prop_Gordon}, and for every
subset $U\subset L_\R$, let 
$$
\cS_{L,\sigma}(U):=\{L\cap(\sigma-v)~|~v\in U\}
\qquad\text{and set}\qquad
\cS_{L,\sigma}:=\cS_{L,\sigma}(L_\R).
$$
There is a natural $L$-module structure on $\cS_{L,\sigma}$;
namely, notice that
$$
(L\cap(\sigma-v))+l=L\cap(\sigma-(v-l))
\qquad
\text{for every $v\in L_\R$ and $l\in L$}
$$
hence the rule $\tau_l:S\mapsto S+l$ defines a bijection of
$\cS_{L,\sigma}$ onto itself, for every $l\in L$, and clearly
$\tau_l\circ\tau_{l'}=\tau_{l+l'}$ for every $l,l'\in L$.
Also, for every $S\in\cS_{L,\sigma}$ define
$$
\Omega(\sigma,S):=\{v\in L_\R~|~L\cap(\sigma-v)=S\}
$$
and denote by $\bar\Omega(\sigma,S)$ the topological
closure of $\Omega(\sigma,S)$ in $L_\R$. For given
$u\in L_\R^\vee$ and $r\in\R$, set
$H_{u,r}:=\{v\in L_\R~|~u(v)\geq r\}$. We shall say that
a subset of $L_\R$ is {\em $\Q$-linearly constructible},
if it lies in the boolean subalgebra of $\cP(L_\R)$
generated by the subsets $H_{u\otimes_\Q\one_\R,r}$,
for $u$ ranging over all the $\Q$-linear forms $L_\Q\to\Q$,
and $r$ ranging over all rational numbers. 

\begin{proposition}\label{prop_linear-constr-part}
With the notation of \eqref{subsec_for-almost-lengths},
the following holds :
\begin{enumerate}
\item
$\cS_{L,\sigma}(U)$ is a finite set, for every bounded subset
$U\subset L_\R$.
\item
$\cS_{L,\sigma}$ is a finitely generated $L$-module.
\item
For every non-empty $S\in\cS_{L,\sigma}$, the subset
$\Omega(\sigma,S)$ is $\Q$-linearly constructible.
\item
Suppose moreover, that $\sigma$ spans $L_\R$. Then, for every
$S\in\cS_{L,\sigma}$, the subset $\Omega(\sigma, S)$ is contained in
the topological closure of its interior (see \eqref{subsec_spann}).
\item
For every $S\in\cS_{L,\sigma}$, and every $v\in\bar\Omega(\sigma,S)$,
we have $S\subset L\cap(\sigma-v)$.
\end{enumerate}
\end{proposition}
\begin{proof} (i): Define $C_\eps$ as in the proof of proposition
\ref{prop_Gordon}; since $U$ is bounded, it is contained in the
union of finitely many subsets of $L_\R$ of the form $C_1+l$,
for $l$ ranging over a finite subset of $L$. On the other hand,
$\tau_l$ induces a bijection
$$
\cS_{L,\sigma}(C_1)\isom\cL_{L,\sigma}(C_1-l)
\qquad
\text{for every $l\in L$}.
$$
Hence, it suffices to check the assertion for $U=C_1$. However,
the proof of proposition \ref{prop_Gordon}(ii) shows that
$L\cap(\sigma-v)$ is generated by $L\cap(C_1-v)$; if $v\in C_1$,
the latter subset is contained in $C':=C_1\cup(-C_1)$, which is
a compact subset of $L_\R$. Therefore $L\cap C'$ is a finite set,
and the claim follows.

(ii): We have already observed that the $L$-module $\cS_{L,\sigma}$
is generated by $\cS_{L,\sigma}(C_1)$, and this is a finite set,
by (i).

(iii): Fix a minimal system $S_1,\dots,S_n$ of generators of
the $L$-module $\cS_{L,\sigma}$ ({\em i.e.} the $S_i$ are
chosen representatives for the orbits of the $L$-action on
$\cS_{L,\sigma}$). After replacing $S_i$ by some translates
$S_i+l$ (for an appropriate $l\in L$) we may also assume that
either $S_i=\emptyset$, or else $0\in S_i$, and notice that
this implies : 
\set\begin{equation}\label{eq_go-to-span}
S_i\subset\La\sigma\Ra\cap L=\sigma_L^\gp
\qquad
\text{for every $i=1,\dots,n$}
\end{equation}
(proposition \ref{prop_Gordon}(i)). Set
$$
A_{ij}:=\{l\in L~|~S_i\subset S_j-a\}
\qquad
\text{for every $i,j\leq n$}
$$
and notice that $A_{ij}$ is a $\sigma_L$-module, for every
$i,j\leq n$.

\begin{claim}\label{cl_A_ij-fg}
If $S_i,S_j\neq\emptyset$, the $\sigma_L$-module $A_{ij}$ is
finitely generated.
\end{claim}
\begin{pfclaim} Fix $l\in L$ such that $\sigma_L+l\subset S_i$.
Next, say that $x_1,\dots,x_t$ is a finite system of generators
for the $\sigma_L$-module $S_j$ (proposition \ref{prop_Gordon}(ii));
by virtue of \eqref{eq_go-to-span}, for every $s=1,\dots,t$, we may
write $x_s=a_s-b_s$ for certain $a_s,b_s\in\sigma_L$. Set
$l':=b_1+\cdots+b_t$, and notice that $S_j\subset\sigma_L-l'$.
Now, if $S_i\subset S_j-a$, we deduce that
$\sigma_L+l\subset\sigma_L-a-l'$, especially $l\in\sigma_L-a-l'$,
{\em i.e.} $a\in\sigma_L-(l+l')$. This shows that $A_{ij}$ is
isomorphic to an ideal of $\sigma_L$, and then the claim follows
from proposition \ref{prop_ideals-in-fg-mon}(ii).
\end{pfclaim}

Now, let $i,j\leq n$ such that $S_i,S_j\neq\emptyset$. Suppose
first that $i\neq j$, and let $A'_{ij}\subset A_{ij}$ be any
finite generating system for the $\sigma_L$-module $A_{ij}$.
From the construction, it is clear that every element of $LS_j$
that contains $S_i$, must contain $S_j-l$, for some $l\in A'_{ij}$.
To deal with the case where $i=j$, we remark, more generally :

\begin{claim}\label{cl_power-bdd}
Let $P$ be any fine and saturated monoid, $M\subset P^\gp$ a
non-empty finitely generated $P$-submodule, and $a\in P^\gp$
an element such that $aM\subset M$. Then $a\in P$.
\end{claim}
\begin{pfclaim} Pick any $m\in M$, and denote by $M'\subset M$
the submodule generated by $(a^km~|~k\in\N)$. According to
proposition \ref{prop_ideals-in-fg-mon}(i), there exists
$N\geq 0$ such that $M'$ is generated by the finite system
$(a^km~|~k=0,\dots,N)$. Especially, $a^{N+1}m\in M'$, and
therefore there exist $x\in P$ and $i\leq N$ such that
$a^{N+1}m=a^imx$ in $M$; it follows that $a^{N+1-i}\in P$,
and finally $a\in P$, since $P$ is saturated.
\end{pfclaim}

From \eqref{eq_go-to-span} we see that $A_{ii}\subset\sigma^\gp_L$,
 if $S_i\neq\emptyset$; combining with claim \ref{cl_power-bdd}, we
deduce that $A_{ii}=\sigma_L$. Moreover, notice as well that
if $S_i=S_i-a$ for some $a\in\sigma^\gp_L$, then both $a$ and
$-a\in A_{ii}$, so that $a\in\sigma^\times_L$. Thus, let $A'_{ii}$
be any set of representatives of
$\fm_\sigma\!\setminus\!\fm_\sigma^2$, where $\fm_\sigma$ denotes
the maximal ideal of $\sigma_L^\sharp$. If $a\in L$, and $S_i-a$
contains strictly $S_i$, then $a$ is a non-invertible element
of $\sigma_L$, and taking into account corollary
\ref{cor_ideals-in-fg-mon}, we see that $A'_{ii}$ is finite,
and there exists $l\in A'_{ii}$ such that $S_i-l\subset S_i-a$.
Next, for every $i\leq n$ such that $S_i\neq\emptyset$, set
$$
\cS^i:=\bigcup_j\{S_j+l~|~l\in A'_{ij}\}
$$
where $j\leq n$ runs over the indices such that $S_j\neq\emptyset$.
Summing up, we conclude that $\cS^i$ is a finite set for every
$i\leq n$ with $S_i\neq\emptyset$, and if an element of
$\cS_{L,\sigma}$ contains strictly $S_i$, then it contains some
element of $\cS^i$. Lastly, in order to prove assertion (iii), we
may assume that $S=S_i$ for some $i\leq n$, and notice that :
\set\begin{equation}\label{eq_describe-these-omegas}
\Omega(\sigma,S_i)=\{v\in L_\R~|~S\subset\sigma-v\}\setminus
\bigcup_{S'\in\cS^i}\{v\in L_\R~|~S'\subset\sigma-v\}.
\end{equation}
Since $S$ is finitely generated (proposition \ref{prop_Gordon}(ii)),
we reduce to showing that, for every $a\in L$, the subset
$\Omega(\sigma,a):=\{v\in L_\R~|~a\in\sigma-v\}$ is
$\Q$-linearly constructible, which follows easily from
proposition \ref{prop_was-part-of-Gordon}(i) and lemma
\ref{lem_ddoble-cone}.

(iv): We remark :

\begin{claim}\label{cl_C_eps-lies}
Let $C_\eps$ be as in the proof of proposition \ref{prop_Gordon}.
Then, for every $S\in\cS_{L,\sigma}$ and every $a\in\Omega(\sigma,S)$
there exists $\eps>0$ such that $a+C_\eps\subset\Omega(\sigma,S)$.
\end{claim}
\begin{pfclaim} Since $\sigma$ is closed in $L_\R$, for every
$a\in L_\R$ and every $b\in L_\R\setminus\Omega(\sigma,a)$
there exists $\eps>0$ such that
$(b+C_\eps)\cap\Omega(\sigma,a)=\emptyset$. Taking into account
\eqref{eq_describe-these-omegas}, the claim follows easily.
\end{pfclaim}

If $\sigma$ span $L_\R$, the subset $C_\eps$ has non-empty
interior $U_\eps$, for every $\eps>0$, and the topological
closure of $U_\eps$ equals $C_\eps$. The assertion is then
an immediate consequence of claim \ref{cl_C_eps-lies}.

(v): The assertion follows easily from proposition
\ref{prop_Gordon}(ii) : the details shall be left to
the reader.
\end{proof}

\sset\subsubsection{}\label{subsec_caffettiera}
Let $L$ be as in \eqref{subsec_for-almost-lengths},
and for all integers $n,m>0$ set
$$
\frac{1}{m}L:=\{v\in L_\Q~|~mv\in L\}
\qquad
\frac{1}{m}L[1/n]:=\bigcup_{k\geq 0}\frac{1}{n^km}L.
$$
For future reference, let us also point out :

\begin{lemma}\label{lem_obama}
With the notation of \eqref{subsec_caffettiera},
let $\Omega\subset L_\R$ be a $\Q$-linearly constructible
subset. Then we have :
\begin{enumerate}
\item
The topological closure of\/ $\Omega$ in $L_\R$ is again
$\Q$-linearly constructible.
\item
There exists an integer $m>0$ such that
$\frac{1}{m}L[1/n]\cap\Omega$ is dense in $\Omega$, for
every $n>1$.
\end{enumerate}
\end{lemma}
\begin{proof} (i): $\Omega$ is a finite union of non-empty
subsets of the form $H_1\cap\cdots\cap H_k$, where each $H_i$
is either of the form $H_{u\otimes\one_\R,r}$ for some non-zero
$\Q$-linear form $u$ of $L_\Q$ and some $r\in\Q$ (and this
is a closed subset of $L_\R$), or else is the complement in
$L_\R$ of a subset of this type (and then its closure is a
half-space $H_{-u\otimes\one_\R,r}$). One verifies that the
closure of $H_1\cap\cdots\cap H_k$ is the intersection of
the closures of $H_1,\dots,H_k$, whence the assertion.

(ii): We may assume that $\Omega=\Omega_1\cap\Omega_2$, where
$\Omega_1$ is a finite intersection of rational hyplerplanes,
and $\Omega_2$ is a finite intersection of open half-spaces
({\em i.e.} of complements of closed half-spaces). Suppose
that $\frac{1}{m}L[1/n]\cap\Omega_1$ is dense in $\Omega_1$;
then clearly $\frac{1}{m}L[1/n]\cap\Omega$ is dense in $\Omega$.
Hence, we may further assume that $\Omega$ is a non-empty
intersection of rational hyperplanes. In this case, $\Omega$
is of the form $V_\R+v_0$, where $v_0\in L_\Q$, and
$V_\R=V\otimes_\Z\R$ for some subgroup $V\subset L$.
Notice that $L[1/n]\cap V_\R$ is dense in $V_\R$ for
every integer $n>1$. Then, any integer $m>0$ such that
$v_0\in\frac{1}{m}L$ will do.
\end{proof}

\subsection{Fine and saturated monoids}\label{sec_fine-satu}
This section presents the further developments of the theory
of fine and saturated monoids. Again, all the monoids in this
section are non-pointed. We begin with a few corollaries
of proposition \ref{prop_Gordon}(i,iii).

\begin{corollary}\label{cor_fragment-Gordon}
Let $M$ be an integral monoid, such that $M^\sharp$ is fine.
We have :
\begin{enumerate}
\item
The inclusion map $M\to M^\sat$ is a finite morphism of monoids.
\item
Especially, if $M$ is fine, any monoid $N$ with
$M\subset N\subset M^\sat$, is fine.
\end{enumerate}
\end{corollary}
\begin{proof} (i): From lemma \ref{lem_exc-satura}(ii) we
deduce that $M^\sat$ is a finitely generated $M$-module if
and only if $(M^\sharp)^\sat$ is a finitely generated
$M^\sharp$-module. Hence, we may replace $M$ by $M^\sharp$,
and assume that $M$ is fine. Pick a surjective group
homomorphism $\phi:\Z^{\oplus n}\to M^\gp$; it is easily
seen that :
$$
\phi^{-1}(M^\sat)=(\phi^{-1}M)^\sat
$$
and clearly it suffices to show that $\phi^{-1}N$ is finitely
generated, hence we may replace $M$ by $\phi^{-1}M$, and assume
throughout that $M^\gp$ is a free abelian group of finite rank. In
this case, proposition \ref{prop_Gordon}(i,iv) already implies
that $M^\sat$ is finitely generated. Let $a_1,\dots,a_k\in M^\sat$
be a finite system of generators, and pick integers $n_1,\dots,n_k>0$
such that $a_i^{n_i}\in M$ for $i=1,\dots,k$. For every $i=1,\dots,k$
let $\Sigma_i:=\{a_i^j~|~j=0,\dots,n_i-1\}$; it is easily seen
$\Sigma_1\cdots\Sigma_k\subset M^\sat$ is a system of generators
for the $M$-module $M^\sat$ (where the product of the sets $\Sigma_i$
is formed in the monoid $\cP(M^\sat)$ of \eqref{subsec_toric}).

(ii) follows from (i), in view of proposition
\ref{prop_ideals-in-fg-mon}(i).
\end{proof}

\begin{corollary}\label{cor_fibres-are-fg}
Let $f:M_1\to M$ and $g:M_2\to M$ be two morphisms of monoids, such
that $M_1$ and $M_2$ are finitely generated, and $M$ is integral.
Then the fibre product $M_1\times_MM_2$ is a finitely generated
monoid, and if $M_1$ and $M_2$ are fine, the same holds for
$M_1\times_MM_2$.
\end{corollary}
\begin{proof} If the monoids $M$, $M_1$ and $M_2$ are integral,
$M_1\times_M M_2$ injects in $M^\gp_1\times_{M^\gp}M_2^\gp$
(lemma \ref{lem_forget-me-not}(iii)), hence it is integral.
To show that the fibre product is finitely generated, choose
surjective morphisms $\N^{\oplus a}\to M_1$ and
$\N^{\oplus b}\to M_2$, for some $a,b\in\N$; by composition we get
maps of monoids $\phi:\N^{\oplus a}\to M$, $\psi:\N^{\oplus b}\to
M$, such that the induced morphism $P:=\N^{\oplus
a}\times_M\N^{\oplus b}\to M_1\times_MM_2$ is surjective. Hence it
suffices to show that $P$ is finitely generated. To this aim, let
$L:=\Ker(\phi^\gp-\psi^\gp:\Z^{\oplus a+b}\to M^\gp)$; for every
$i=1,\dots,a+b$, denote also by $\pi_i:\Z^{\oplus a+b}\to\Z$ the
projection onto the $i$-th direct summand. The system
$\{\pi_i~|~i=1,\dots,a+b\}$ generates a rational convex polyhedral
cone $\sigma\subset L^\vee\otimes_\Z\R$, and one verifies easily
that $P=L\cap\sigma^\vee$, so the assertion follows from propositions
\ref{prop_was-part-of-Gordon}(i) and \ref{prop_Gordon}(i).
\end{proof}

\begin{corollary}\label{cor_no-fibres-here}
Let $(\Gamma,+,0)$ be an integral monoid, $M$ a finitely generated\/
$\Gamma$-graded monoid. Then $M_0$ is a finitely generated monoid,
and $M_\gamma$ is a finitely generated $M_0$-module, for every
$\gamma\in\Gamma$.
\end{corollary}
\begin{proof} We have $M_0=M\times_\Gamma\{0\}$, hence $M_0$
is finitely generated, by corollary \ref{cor_fibres-are-fg}.
The given element $\gamma\in\Gamma$ determines a
unique morphism of monoids $\N\to\Gamma$ such that $1\mapsto\gamma$.
Let $p_1:M':=M\times_\Gamma\N\to\N$ and $p_2:M'\to M$ be the two
natural projections; by lemma \ref{lem_forget-me-not}(iii), we
have $M_\gamma=p_2(p_1^{-1}(1))$. In light of corollary
\ref{cor_fibres-are-fg}, $M'$ is still finitely generated,
hence we are reduced to the case where $\Gamma=\N$ and $\gamma=1$.
In this case, pick a finite set of generators $S$ for $M$. One
checks easily that $M_1\cap S$ generates the $M_0$-module $M_1$.
\end{proof}

\begin{corollary}\label{cor_satu-gives-flat}
Let $M$ be an integral monoid, such that $M^\sharp$ is fine,
and $\phi:M\to N$ a saturated morphism of monoids. Then $\phi$
is flat.
\end{corollary}
\begin{proof} In view of corollary \ref{cor_fragment-Gordon}(i) and
theorem \ref{th_Gruson}, it suffices to show that the $M^\sat$-module
$M^\sat\otimes_MN$ is flat. Hence we may replace $M$ by $M^\sat$,
and assume that $M$ is saturated.
Let $I\subset M$ be any ideal, and define $\sR(M,I)$ as in
the proof of theorem \ref{th_flat-crit-for-mnds}; by assumption,
$\sR(M,I)^\sat\otimes_MN$ is a saturated -- especially, integral --
monoid, {\em i.e.} the natural map
$$
\sR(M,I)\otimes_MN\to\sR(M,I)^\gp\otimes_{M^\gp}N^\gp
$$
is injective. The latter factors through the morphism $j\otimes_MN$,
where $j:\sR(M,I)\to M\times\N$ is the obvious inclusion.
In light of example \ref{ex_Rees-satura}(i), we deduce that
the induced map $I^\sat\otimes_MN\to N$ is injective.
Now, if $I$ is a prime ideal, then $I^\sat=I$, hence the
contention follows from corollary \ref{cor_first-corollo}(ii).
\end{proof}

The following corollary generalizes lemma \ref{lem_decomp-sats}.

\begin{corollary}\label{cor_general-deco}
Let $f:M\to N$ be a local, flat and saturated morphism of fine
monoids, with $M$ sharp. Then there exists an isomorphism of monoids
$$
g:N\isom N^\sharp\times N^\times
$$
that fits into a commutative diagram
$$
\xymatrix{ M \ar[r]^-{f^\sharp} \ar[d]_f & N^\sharp \ar[d] \\
           N \ar[r]^-g & N^\sharp\times N^\times
}$$
whose right vertical arrow is the natural inclusion map.
\end{corollary}
\begin{proof} From lemma \ref{lem_persist-integr}(ii), we know
that $f$ is exact, and since $M$ is sharp, we easily deduce
that $f(M)^\gp\cap N^\times=\{1\}$. Hence, the induced group
homomorphism $M^\gp\oplus N^\times\to N^\gp$ is injective. On
the other hand, since $f$ is flat, local and saturated, the same
holds for $f^\sharp:M\to N^\sharp$ (lemma \ref{lem_little}(iii)
and corollary \ref{cor_satu-gives-flat}); then corollary
\ref{cor_persist-integr}(ii) says that the cokernel of the
induced group homomorphism
$M^\gp\to (N^\sharp)^\gp=N^\gp/N^\times$ is a free abelian
group $G$ (of finite rank). Summing up, we obtain an isomorphism
of abelian groups :
$$
h:M^\gp\oplus N^\times\oplus G\isom N^\gp
$$
extending the map $f^\gp$. Set $N_0:=N\cap h(M^\gp\oplus G)$;
it follows easily that the natural map $N_0\times N^\times\to N$
is an isomorphism; especially, the projection $N\to N^\sharp$
maps $N_0$ isomorphically onto $N^\sharp$, and the contention
follows.
\end{proof}

\sset\subsubsection{}\label{subsec_Gordon-more} Let $(M,\cdot)$ be a
fine (non-pointed) monoid, so that $M^\gp$ is a finitely generated
abelian group. We set $M^\gp_\R:=\log M^\gp\otimes_\Z\R$, and we let
$M_\R$ be the convex polyhedral cone generated by the image of $\log
M$. Then $(M_\R,+)$ is a monoid, and we have a natural morphism of
monoids
$$
\phi:\log M\to(M_\R,+).
$$

\begin{proposition}\label{prop_Gordie}
With the notation of \eqref{subsec_Gordon-more}, we have :
\begin{enumerate}
\item
The rule : $F\mapsto F_\R$ establishes a bijection between the set
of faces of $M$ and the set of faces of $M_\R$.
\item
The map $\phi$ induces a bijection
$$
\phi^*:\Spec\,M_\R\to\Spec\,M.
$$
\end{enumerate}
\end{proposition}
\begin{proof} Clearly, we may assume that $M\neq\{1\}$.
Let $\fp\subset M$ be a prime ideal; we denote by $\fp_\R$ the
ideal of $(M_\R,+)$ generated by all elements of the form
$r\cdot\phi(x)$, where $r$ is any strictly positive real number, and
$x$ is any element of $\fp$. We also denote by
$(M\!\setminus\!\fp)_\R$ the convex cone of $M^\gp_\R$ generated by
the image of $M\!\setminus\!\fp$.

\begin{claim}\label{cl_disj-cones}
$M_\R$ is the disjoint union of $(M\!\setminus\!\fp)_\R$ and
$\fp_\R$.
\end{claim}
\begin{pfclaim}
To begin with, we show that $M_\R=(M\!\setminus\!\fp)_\R\cup\fp_\R$.
Indeed, let $x\in M_\R$; then we may write $x=\sum_{i=1}^hm_i\otimes
a_i$ for certain $a_1,\dots,a_h\in\R_+$ and $m_1,\dots,m_h\in M$. We
may assume that $a_1,\dots,a_k\in\fp$ and $a_{k+1},\dots,a_h\in
M\!\setminus\!\fp$. Now, if $k=0$ we have
$x\in(M\!\setminus\!\fp)_\R$, and otherwise $x\in\fp_\R$, which
shows the assertion.

It remains to show that
$(M\!\setminus\!\fp)_\R\cap\fp_\R=\emptyset$. To this aim, suppose
by way of contradiction, that this intersection contains an element
$x$; this means that we have finite subsets $S_0\subset
M\!\setminus\!\fp$ and $S_1\subset M$ such that
$S_1\cap\fp\neq\emptyset$, and an identity of the form :
\set\begin{equation}\label{eq_win-edt}
x=\sum_{\sigma\in S_0}\sigma\otimes a_\sigma=
\sum_{\sigma\in S_1}\sigma\otimes b_\sigma
\end{equation}
where $a_\sigma>0$ for every $\sigma\in S_0$ and $b_\sigma>0$ for
every $\sigma\in S_1$. For every $\sigma\in S_0$, choose a rational
number $a'_\sigma\geq a_\sigma$; after adding the summand
$\sum_{\sigma\in S_0}\sigma\otimes (a'_\sigma-a_\sigma)$ to both
sides of \eqref{eq_win-edt}, we may assume that $a_\sigma\in\Q_+$
for every $\sigma\in S_0$. Let $N\subset M$ be the submonoid
generated by $S_1$; it follows that $x\in N_\R\cap M_\Q=N_\Q$
(proposition \ref{prop_Gordon}(iii)), hence we may assume that
all the coefficients $a_\sigma$ and $b_\sigma$ are rational
and strictly positive (see remark \ref{cl_like-this-x}).
We may further multiply both sides of \eqref{eq_win-edt} by
a large integer, to obtain that these coefficients are actually
integers. Then, up to further multiplication by some integer,
the identity of \eqref{eq_win-edt} lifts to an identity between
elements of $\log M$, of the form :
$\sum_{\sigma\in S_0}a_\sigma\cdot\sigma=\sum_{\sigma\in
S_1}b_\sigma\cdot\sigma$. The latter is absurd, since
$S_1\cap\fp\neq\emptyset$ and $S_0\cap\fp=\emptyset$.
\end{pfclaim}

Claim \ref{cl_disj-cones} implies that $\fp_\R$ is a prime ideal of
$M_\R$, and clearly $\fp\subset\phi^*(\fp_\R)$. Since we have as
well $M\!\setminus\!\fp\subset\phi^{-1}(M\!\setminus\!\fp)_\R$, we
deduce that $\fp=\phi^*(\fp_\R)$. Hence the rule $\fp\mapsto\fp_\R$
yields a right inverse $\phi_*:\Spec\,M\to\Spec\,M_\R$ for the
natural map $\phi^*$. To show that $\phi_*$ is also a left inverse,
let $\fq\subset M_\R$ be a prime ideal; by lemma \ref{lem_same-faces}
and proposition \ref{prop_was-part-of-Gordon}(i), the face
$M_\R\!\setminus\!\fq$ is of the form $M_\R\cap\Ker\,u$,
for some $u\in M^\vee_\R\cap(\log M^\gp)^\vee$. Then it is
easily seen that $M_\R\!\setminus\!\fq$ is the convex cone
generated by $\phi(M)\cap\Ker\,u$, in other words,
$M_\R\setminus\fq=\phi^{-1}(\Ker\,u)_\R$.
Again by claim \ref{cl_disj-cones}, it follows that
$\fq=(M\!\setminus\!\phi^{-1}(\Ker\,u))_\R=(\phi^*\fq)_\R$, as
stated. The argument also shows that every face of $M_\R$ is
of the from $(M\!\setminus\!\fp)_\R$ for a unique prime ideal
$\fp$, which settles assertion (i).
\end{proof}

\begin{corollary}\label{cor_consequent}
Let $M$ be a fine monoid. We have :
\begin{enumerate}
\item
$\dim M=\rk_\Z(M^\gp/M^\times)$.
\item
$\dim(M\!\setminus\!\fp)+\hgt\,\fp=\dim M$ \ \ for every
$\fp\in\Spec\,M$.
\item
If $M\neq\{1\}$ is sharp (see \eqref{subsec_regular-mon}), there
exists a local morphism $M\to\N$.
\item
If $M$ is sharp and $M^\gp$ is a torsion-free abelian group of
rank $r$, there exists an injective morphism of monoids
$M\to\N^{\oplus r}$.
\end{enumerate}
\end{corollary}
\begin{proof}
(i): By proposition \ref{prop_Gordie}, the dimension of $M$ can be
computed as the length of the longest chain $F_0\subset
F_1\subset\cdots\subset F_d$ of strict inclusions of faces of
$M_\R$. On the other hand, given such a maximal chain, denote by
$r_i$ the dimension of the $\R$-vector space spanned by $F_i$; in
view proposition \ref{prop_refer-to-Ful}(ii),(iii), it is easily
seen that $r_{i+1}-r_i=1$ for every $i=0,\dots,d-1$. Since
$M_\R\cap(-M_\R)$ is the minimal face of $M_\R$, we deduce that
$$
\dim M=\dim_\R M^\gp_\R-\dim_\R M_\R\cap(-M_\R).
$$
Clearly $\dim_\R M^\gp_\R=\rk_\Z M^\gp$; moreover, by proposition
\ref{prop_Gordie}, the face $M_\R\cap(-M_\R)$ is spanned by the
image of the face $M^\times$ of $M$. whence the assertion.

(ii) is similar : again proposition
\ref{prop_refer-to-Ful}(ii),(iii) implies that, every face $F$ of
$M_\R$ fits into a maximal strictly ascending chain of faces of
$M_\R$, and the length of any such maximal chain is $\dim M$, by
(i).

(iii): Notice that $\rk_\Z M^\gp>0$, by (i). By proposition
\ref{prop_Gordie}(i), $M_\R$ is strictly convex, therefore,
by proposition \ref{prop_was-part-of-Gordon}(i), we may find
a non-zero linear map $\phi:M^\gp\otimes_\Z\Q\to\Q$, such that
$M_\R\cap\Ker\,\phi\otimes_\Q\R=\{0\}$ and $\phi(M)\subset\Q_+$.
A suitable positive integer of $\phi$ will do.

(iv): Under the stated assumption, we may regard $M$ as a submonoid
of $M_\R$, and the latter contains no non-zero linear subspaces. By
corollary \ref{cor_strongly} and proposition
\ref{prop_was-part-of-Gordon}(i), we may then find $r$ linearly
independent forms $u_1,\dots,u_r:M^\gp\otimes_\Z\Q\to\Q$ which are
positive on $M$. It follows that
$u_1\otimes_\Q\R,\dots,u_r\otimes_\Q\R$ generate a polyhedral cone
$\sigma^\vee\subset M^\vee_\R$, so its dual cone
$\sigma\subset M^\gp_\R$ contains $M_\R$. By construction, $\sigma$
admits precisely $r$ extremal rays, say the rays generated by the
vectors $v_1,\dots,v_r$, which we can pick in $M^\gp_\Q$, in which
case they form a basis of the latter $\Q$-vector space. Now, every
$x\in M_\R$ can be written uniquely in the form
$x=\sum_{i=1}^ra_iv_i$ for certain $a_1,\dots,a_r\in\Q_+$; since $M$
is finitely generated, we may find an integer $N>0$ independent of
$x$, such that $Na_i\in\N$ for every $i=1,\dots,r$. In other words,
$M$ is contained in the monoid generated by
$N^{-1}v_1,\dots,N^{-1}v_r$; the latter is isomorphic to
$\N^{\oplus r}$.
\end{proof}

\sset\subsubsection{}\label{subsec_dual-of-mon}
For any monoid $M$, the {\em dual\/} of $M$ is the monoid
$$
M^\vee:=\Hom_\Mnd(M,\N)
$$
(see \eqref{subsec_toric}). As usual, there is a natural
morphism
$$
M\to M^{\vee\vee}
\quad : \quad
m\mapsto(\phi\mapsto\phi(m))
\quad
\text{for every $m\in M$ and $\phi\in M^\vee$}.
$$
We say that $M$ is {\em reflexive\/}, if this morphism is an
isomorphism.

\begin{proposition}\label{prop_reflex-dual}
Let $M$ be a monoid. We have :
\begin{enumerate}
\item
$M^\vee$ is integral, saturated and sharp.
\item
If $M$ is finitely generated, $M^\vee$ is fine, and we have a
natural identification :
$$
(M^\vee)_\R\isom(M_\R)^\vee
$$
(where $(M_\R)^\vee$ is defined as in \eqref{subsec_convex-cone}).
Moreover,\ \  $\dim M=\dim M^\vee$.
\item
If $M$ is finitely generated and sharp, we have a natural
identification :
$$
(M^\vee)^\gp\isom(M^\gp)^\vee.
$$
\item
If $M$ is fine, sharp and saturated, then $M$ is reflexive.
\end{enumerate}
\end{proposition}
\begin{proof}(i): It is easily seen that the natural
group homomorphism
\set\begin{equation}\label{eq_dual-and-gp}
(M^\vee)^\gp\to(M^\gp)^\vee:=\Hom_\Z(M^\gp,\Z)
\end{equation}
is injective. Now, say that $\phi\in(M^\vee)^\gp$ and
$N\phi\in M^\vee$ for some $N\in\N$; we may view $\phi$ as
group homomorphism $\phi:M^\gp\to\Z$, and the assumption
implies that $\phi(M)\subset\Z\cap\Q_+=\N$, whence the
contention.

(ii): Indeed, let $x_1,\dots,x_n$ be a system of generators
of $M$. Define a group homomorphism
$f:(M^\gp)^\vee\to\Z^{\oplus n}$ by the rule :
$\phi\mapsto(\phi(x_1),\dots,\phi(x_n))$ for every
$f:M^\gp\to\Z$. Then $M^\vee=\phi^{-1}(\N^{\oplus n})$, and since
$(M^\gp)^\vee$ is fine, corollary \ref{cor_fibres-are-fg} implies
that $M^\vee$ is fine as well. Next, the injectivity of
\eqref{eq_dual-and-gp} implies especially that $(M^\vee)^\gp$
is torsion-free, hence $\eqref{eq_dual-and-gp}\otimes_\Z\R$
is still injective; its restriction to $(M^\vee)_\R$ factors
therefore through an injective map $f:(M^\vee)_\R\to(M_\R)^\vee$.
The latter map is determined by the image of $M^\vee$, and
by inspecting the definitions, we see that
$f(\phi):=\phi^\gp\otimes 1$ for every $\phi\in M^\vee$.
To prove that $f$ is an isomorphism, it suffices to show that
it has dense image. However, say that $\phi\in(M_\R)^\vee$;
then $\phi:M^\gp\to\R$ is a group homomorphism such that
$\phi(M)\subset\R_+$. Since $M$ is finitely generated, in
any neighborhood of $\phi$ in $M^\gp_\R$ we may find some
$\phi':M^\gp\to\Q_+$, and then $N\phi'\in M^\vee$ for
some integer $N\in\N$ large enough. It follows that $\phi'$
is in the image of $f$, whence the contention.

The stated equality follows from the chain of identities :
$$
\dim M=\dim M_\R=\dim (M_\R)^\vee=\dim (M^\vee)_\R=\dim M^\vee
$$
where the first and the last follow from proposition
\ref{prop_Gordie}(ii), and the second follows from corollary
\ref{cor_spanning-cone}(ii).

(iii): Let us show first that, under these assumptions,
$\eqref{eq_dual-and-gp}\otimes_\Z\R$ is an isomorphism.
Indeed, if $M$ is sharp, $(M_\R)^\vee$ spans $(M^\gp_\R)^\vee$
(corollary \ref{cor_strongly} and proposition \ref{prop_Gordie}(i));
then the assertion follows from (ii).
We deduce that $(M^\vee)^\gp$ and $(M^\gp)^\vee$ are free
abelian groups of the same rank, hence we may find a basis
$\phi_1,\dots,\phi_r$ of $(M^\vee)^\gp$ (resp.
$\psi_1,\dots,\psi_r$ of $(M^\gp)^\vee$), and positive integers
$N_1,\dots,N_r$ such that \eqref{eq_dual-and-gp} is given by
the rule : $\phi_i\mapsto N_i\psi_i$ for every $i=1,\dots,r$.
But then necessarily we have $N_i=1$ for every $i\leq r$, and
(iii) follows.

(iv): It is easily seen that
$M^\vee=(M_\R)^\vee\cap(M^\gp)^\vee$ (notation of
\eqref{subsec_Gordon-more}). After dualizing again we find :
$M^{\vee\vee}=((M^\vee)_\R)^\vee\cap(M^{\vee\gp})^\vee$.
From (ii) we deduce that $((M^\vee)_\R)^\vee=(M_\R)^{\vee\vee}=M_\R$
(lemma \ref{lem_ddoble-cone}), and from (iii) we get :
$(M^{\vee\gp})^\vee=(M^\gp)^{\vee\vee}=M^\gp$.
Hence $M^{\vee\vee}=M_\R\cap M^\gp=M$ (proposition
\ref{prop_Gordon}(iii)).
\end{proof}

\begin{remark}\label{rem_reflex-dual}
(i)\ \ Let $M$ be a sharp and fine monoid. Proposition
\ref{prop_reflex-dual}(iii) implies that the natural map
$$
\Hom_\Mnd(M,\Q_+)^\gp\to\Hom_\Mnd(M,\Q)
$$
is an isomorphism. Indeed, it is easily seen that this map is
injective. For the surjectivity, one uses the identification
$\Hom_\Mnd(M,\Q)\isom(M^\vee)^\gp\otimes_\Z\Q$, which follows from
{\em loc.cit.} (Details left to the reader.)

(ii)\ \ For $i=1,2$, let $N_i\to N$ be two morphisms of monoids. By
general nonsense, we have a natural isomorphism :
$$
(N_1\otimes_N N_2)^\vee\isom N^\vee_1\times_{N^\vee}N_2^\vee.
$$

(iii)\ \ If $f_i:M_i\to M$ ($i=1,2$) are morphisms of fine,
saturated and sharp monoids, there exists a natural surjection :
\set\begin{equation}\label{eq_double_dual_amalg}
M^\vee_1\otimes_{M^\vee}M^\vee_2\to(M_1\times_MM_2)^\vee
\end{equation}
whose kernel is the subgroup of invertible elements.
Indeed, set $P:=M^\vee_1\otimes_{M^\vee}M_2^\vee$; in view of (ii)
and proposition \ref{prop_reflex-dual}(iv), we have a natural
identification $P^{\vee\vee}\isom(M_1\times_MM_2)^\vee$, and the
sought map is its composition with the double duality map
$P\to P^{\vee\vee}$.
Moreover, clearly $P$ is finitely generated, and it is also
integral and saturated, since saturation commutes with colimits.
Hence -- again by proposition \ref{prop_reflex-dual}(iv) -- the
double duality map induces an isomorphism $P/P^\times\isom P^{\vee\vee}$.

(iv)\ \ In the situation of (iii), if $f_i:M_i\to M$ ($i=1,2$) are
epimorphisms, then \eqref{eq_double_dual_amalg} is an isomorphism.
Indeed, in this case the dual morphisms
$f^\vee_i:M^\vee\to M^\vee_i$ are injective, so that $P$ is sharp
(lemma \ref{lem_localize-sharp}), whence the claim.
\end{remark}

\begin{theorem}\label{th_structure-of-satu}
Let $M$ be a saturated monoid, such that $M^\sharp$ is fine.
We have :
\begin{enumerate}
\item
$M=\displaystyle{\bigcap_{\hgt\,\fp=1}M_\fp}$\ \  (where the
intersection runs over the prime ideals of $M$ of height one).
\item
If moreover, $\dim M=1$, then there is an isomorphism of monoids :
$$
M^\times\times\N\isom M.
$$
\item
Suppose that $M^\gp$ is a torsion-free abelian group, and let $R$ be
any normal domain. Then the group algebra $R[M]$ is a normal domain
as well.
\end{enumerate}
\end{theorem}
\begin{proof}(i): Pick a decomposition $M=M^\sharp\times M^\times$
as in lemma \ref{lem_decomp-sats}, and notice that $M^\sharp$ is
fine, sharp and saturated. The prime ideals of $M$ are of the form
$\fp=\fp_0\times M^\times$, where $\fp_0$ is a prime ideal of $M^\sharp$.
Then it is easily seen that $M_\fp=M^\sharp_{\fp_0}\times M^\times$.
Therefore, the sought assertion holds for $M$ if and only if it holds
for $M^\sharp$, and therefore we may replace $M$ by $M^\sharp$, which
reduces to the case where $M$ is sharp, hence the natural morphism
$\phi:\log M\to M_\R$ is injective. In such situation, we have
$M=M_\R\cap M^\gp$ and $M_\fp=M_{\fp,\R}\cap M^\gp$ for every prime
ideal $\fp\subset M$ (proposition \ref{prop_Gordon}(iii) and lemma
\ref{lem_exc-satura}(i)). Thus, we are reduced to showing that
$$
M_\R=\bigcap_{\hgt\,\fp=1}M_{\fp,\R}.
$$
However, set $\tau:=(M\!\setminus\!\fp)_\R$; by inspecting the
definitions, one sees that $M_{\fp,\R}=M_\R+(-\tau)$, and
proposition \ref{prop_Gordie} shows that $\tau$ is a facet of
$M_\R$, hence $M_{\fp,\R}$ is the half-space denoted $H_\tau$ in
\eqref{subsec_spann}. Then the assertion is a rephrasing of
proposition \ref{prop_spanning-cone}(ii).

(ii): Arguing as in the proof of (i), we may reduce again to the
case where $M$ is sharp, in which case $M=M_\R\cap M^\gp$. The
foregoing shows that, in case $\dim M=1$, the cone
$M_\R$ is a half-space, whose boundary hyperplane is the only
non-trivial face $\sigma$ of $M_\R$. However, $\sigma$ is generated
by the image of the unique non-trivial face of $M$, {\em i.e.}
by $M^\times=\{1\}$ (proposition \ref{prop_Gordie}(i)), hence
$\sigma=\{0\}$, so $M_\R$ is a half-line. Now, let $u:M^\gp_\R\to\R$
be a non-zero linear form, such that $u(M)\geq 0$, and $x_1,\dots,x_n$
a system of non-zero generators for $M$; say that $u(x_1)$ is the least
of the values $u(x_i)$, for $i=1,\dots,n$. Since $M$ is saturated,
it follows easily that every value $u(x_i)$ is an integer multiple
of $u(x_1)$ (proposition \ref{prop_Gordon}(iii)), and then $x_1$ is
a generator for $M$, so $M\simeq\N$.

(iii): To begin with, $R[M]\subset R[M^\gp]$, and since $M^\gp$ is
torsion-free, it is clear that $[M^\gp]$ is a domain, hence the same
holds for $R[M]$. Furthermore, from (i) we derive :
$R[M]=\bigcap_{\hgt\fp=1}R[M_\fp]$, hence it suffices to show
that $R[M_\fp]$ is normal whenever $\fp$ has height one. However, we
have $R[M_\fp]\simeq R[M^\times_\fp]\otimes_RR[\N]$ in light of
(ii), and since $M^\times_\fp$ is torsion free, it is
a filtered colimit of a family of free abelian groups of finite
rank, so everything is clear.
\end{proof}

\begin{example}\label{ex_satu-dim-two}
Let $M$ be a fine, sharp and saturated monoid of dimension $2$. 

(i)\ \
By corollary \ref{cor_consequent}(i) and example \ref{ex_cone-dim-two},
we see that $M$ admits exactly two facets, which are fine saturated
monoids of dimension one; by theorem \ref{th_structure-of-satu}(ii)
each of these facets is generated by an element, say $e_i$ (for $i=1,2$).
From proposition \ref{prop_Gordon}(iii) it follows that
$\Q_+e_1\oplus\Q_+e_2=M_\Q$. Especially, we may find an integer
$N>0$ large enough, such that :
$$
\N e_1\oplus\N e_2\subset
M\subset\N\frac{e_1}{N}\oplus\N\frac{e_2}{N}.
$$

(ii)\ \
Moreover, clearly $e_1$ and $e_2$ are {\em unimodular\/}
elements of $M^\gp$ ({\em i.e.} they generate direct summands
of the latter free abelian group of rank $2$). We may then
find a basis $f_1,f_2$ of $M^\gp$ with $e_1=f_1$, and
$e_2=af_1+bf_2$, where $a,b\in\Z$ and $(a,b)=1$. After
replacing $f_2$ by some element of the form $cf_2+df_1$
with $c\in\{1,-1\}$ and $d\in\Z$, we may assume that
$b>0$ and $0\leq a<b$. Clearly, such a normalized pair
$(a,b)$ determines the isomorphism class of $M$, since
$M_\R$ is the strictly convex cone of $M^\gp_\R$ whose
extremal rays are generated by $e_1$ and $e_2$, and
$M=M^\gp\cap M_\R$.

(iii)\ \
More precisely, suppose that $M'$ is another fine, sharp and
saturated monoid of dimension $2$, and $\phi:M\to M'$ an
isomorphism. Pick a basis $f_1',f_2'$ of $M'{}^\gp$ and a
normalized pair $(a',b')$ as in (ii), such that $e'_1:=f_1$
and $e'_2:=a'f'_1+b'f'_2$ generate the two facets of $M'$.
Clearly, $\phi$ must send a facet of $M$ onto a facet of
$M'$; we distinguish two possibilities :
\begin{itemize}
\item
either $\phi(e_1)=e'_1$ and $\phi(e_2)=e'_2$, in which
case we get $\phi(f_2)=b^{-1}(a'-a)f'_1+b^{-1}b'f'_2$;
especially, $b',a-a'\in b\Z$. By considering $\phi^{-1}$,
we get symmetrically that $b\in b'\Z$, so $b=b'$ and
therefore $(a',b)=1=(a,b)$ and $0\leq a,a'<b$, whence
$a=a'$
\item
or else $\phi(e_1)=e'_2$ and $\phi(e_2)=e'_1$, in which
case we get $\phi(f_2)=b^{-1}(1-aa')f_1+b^{-1}b'af_2$.
It follows again that $b'\in b\Z$, so $b=b'$, arguing as
in the previous case. Moreover, $0\leq a'<b$, and
$1-aa'\in b\Z$. In other words, the class of $a'$ in
the group $(\Z/b\Z)^\times$ is the inverse of the class
of $a$.
\end{itemize}
Conversely, it is easily seen that if $M'$ is as above,
$f'_1,f'_2$ is a basis of $M'{}^\gp$, and the two facets
of $M'$ are generated by $f'_1$ and $a'f'_1+b'f'_2$, for
a pair $(a',b')$ normalized as in (ii), and such that
$aa'\equiv 1{\pmod b}$, then there exists an isomorphism
$M\isom M'$ of monoids (details left to the reader).
Hence, set $(\Z/b\Z)^\dagger:=(\Z/b\Z)^\times/\!\!\sim$, where
$\sim$ denotes the smallest equivalence relation such that
$[a]\sim[a]^{-1}$ for every $[a]\in(\Z/b\Z)^\times$. 
We conclude that there exists a natural bijection between
the set of isomorphism classes of fine, sharp and saturated
monoids of dimension $2$, and the set of pairs $(b,[a])$,
where $b>0$ is an integer, and $[a]\in(\Z/b\Z)^\dagger$.
\end{example}

\sset\subsubsection{}\label{subsec_fract-ideals}
Let $P$ be an integral monoid. A {\em fractional ideal\/} of $P$
is a $P$-submodule $I\subset P^\gp$ such that $I\neq\emptyset$
and $x\cdot I\subset P$ for some $x\in P$. Clearly the union and
the intersection of finitely many fractional ideals, are again
fractional ideals. We may also define the product of two fractional
ideals $I_1,I_2\subset P^\gp$ : namely, the subset
$$
I_1I_2:=\{xy~|~x\in I_1, y\in I_2\}\subset P^\gp
$$
which is again a fractional ideal, by an easy inspection. If $I$ is a
fractional ideal of $P$, we say that $I$ is {\em finitely generated},
if it is such, when regarded as a $P$-module. For any two fractional
ideals $I_1,I_2$, we let
$$
(I_1:I_2):=\{x\in P^\gp~|~x\cdot I_2\subset I_1\}.
$$
It is easily seen that $(I_1:I_2)$ is a fractional ideal of $P$
(if $x\in I_2$ and $yI_1\subset P$, then clearly $xy(I_1:I_2)\subset P$).
We set
$$
I^{-1}:=(P:I)
\qquad\text{and}\qquad
I^*:=(I^{-1})^{-1}
\qquad
\text{for every fractional ideal $I\subset P^\gp$}.
$$
Clearly $J^{-1}\subset I^{-1}$, whenever $I\subset J$, and
$I\subset I^*$ for all fractional ideals $I$, $J$. We say
that $I$ is {\em reflexive\/} if $I=I^*$.  We remark that
$I^{-1}$ is reflexive, for every fractional ideal
$I\subset P^\gp$. Indeed, we have $I^{-1}\subset(I^{-1})^*$,
and on the other hand $(I^{-1})^*=(I^*)^{-1}\subset I^{-1}$.
It follows that $I^*$ is reflexive, for every fractional
ideal $I$. Moreover, $I^*\subset J^*$, whenever $I\subset J$;
especially, $I^*$ is the smallest reflexive fractional
ideal containing $I$. Notice furthermore, that
$aI^{-1}=(a^{-1}I)^{-1}$ for every $a\in P^\gp$; therefore,
$aI^*=(aI)^*$, for every fractional ideal $I$ and $a\in P^\gp$.

\begin{lemma}\label{lem_associat}
Let $P$ be any integral monoid, $I,J\subset P^\gp$ two
fractional ideals. Then :
\begin{enumerate}
\item
$(IJ)^*=(I^*J^*)^*$.
\item
$I^*$ is the intersection of the invertible fractional
ideals of $P$ that contain $I$ (see definition
{\em\ref{def_coh-idea-log}(iv)}).
\end{enumerate}
\end{lemma}
\begin{proof}(i): Since $IJ\subset I^*J^*$, we have
$(IJ)^*\subset(I^*J^*)^*$. To show the converse inclusion,
it suffices to check that $I^*J^*\subset(IJ)^*$, since
$(IJ)^*$ is reflexive, and $(I^*J^*)^*$ is the smallest
reflexive fractional ideal containing $I^*J^*$. Now, let
$a\in I$ be eny element; we get $aJ^*=(aJ)^*\subset(IJ)^*$,
so $IJ^*\subset(IJ)^*$, and therefore $(IJ^*)^*\subset(IJ)^*$.
Lastly, let $b\in J^*$ be any element; we get
$bI^*=(bI)^*\subset(IJ^*)^*$, so $I^*J^*\subset(IJ^*)^*$,
whence the lemma.

(ii): It suffices to unwind the definitions. Indeed,
$a\in P^\gp$ lies in $I^*$ if and only if $aI^{-1}\subset P$,
if and only if $ab\in P$, for every $b\in P^\gp$ such that
$bI\subset P$. In other words, $a\in I^*$ if and only if
$a\in b^{-1}P$ for every $b\in P^\gp$ such that $I\subset b^{-1}P$,
which is the contention.
\end{proof}

\sset\subsubsection{}\label{subsec_Div}
Let $P$ be any integral monoid. We denote by
$\Div(P)$ the set of all reflexive fractional ideals of $P$.
We define a composition law on $\Div(P)$ by the rule :
$$
I\odot J:=(IJ)^*
\qquad
\text{for every $I,J\in\Div(P)$}.
$$
It follows easily from lemma \ref{lem_associat}(i) that $\odot$
is an associative law; indeed we may compute :
$$
(I\odot J)\odot K=((IJ)^*K)^*=(IJK)^*=(I(JK)^*)^*
=I\odot(J\odot K)
$$
for every $I,J,K\in\Div(P)$. Clearly $I\odot J=J\odot I$ and
$P\odot I=I$, for every $I,J\in\Div(P)$, so $(\Div(P),\odot)$
is a commutative monoid. Notice as well that if $I\subset P$,
then also $I^*\subset P$ (lemma \ref{lem_associat}(ii)), so
the subset of all reflexive fractional ideals contained in $P$
is a submonoid $\Div_+(P)\subset\Div(P)$.

\begin{example}\label{ex_integral-doms}
Let $A$ be an integral domain, and $K$ the field of fractions
of $A$. Classically, one defines the notion of {\em fractional
ideal} : see {\em e.g.} \cite[p.80]{Mat} and
\cite[Ch.VII, \S1, no.1]{BouAC}, but notice that the definitions
in these two references differ slightly, as the zero ideal is
a fractional ideal according to the latter, but not according
to the former. Additionally, one has a notion of {\em reflexive
fractional ideal\/} of $A$, which are also called {\em divisorial
fractional ideals} in \cite[Ch.VII, \S1, no.1]{BouAC}. In our
terminology, these are understood as follows. Set $A':=A\cap K^\times$,
and notice that the monoid $(A,\cdot)$ is naturally isomorphic
to the integral pointed monoid $A'_\circ$. Then a fractional ideal
of $A$ is an $A'_\circ$-submodule of $K^\times_\circ=K$ of the form
$I_\circ$, where $I\subset K^\times$ is a fractional ideal of $A'$.
Likewise one may define the reflexive fractional ideals of $A$.
The set $\Div(A)$ of all reflexive fractional ideals of $A$ is
then endowed with the unique monoid structure, such that the map
$\Div(A')\to\Div(A)$ given by the rule $I\mapsto I_\circ$ is an
isomorphism of monoids.
\end{example}

\begin{lemma}\label{lem_rflx_and_quot}
Let $P$ be an integral monoid, and $G\subset P^\times$ a subgroup.
We have :
\begin{enumerate}
\item
The rule $I\mapsto I/G$ induces a bijection from the set of
fractional ideals of $P$ to the set of fractional ideals of $P/G$.
\item
A fractional ideal $I$ of $P$ is reflexive if and only if the same
holds for $I/G$.
\item
The rule $I\mapsto I/G$ defines an isomorphism of monoids
$$
\Div(P)\isom\Div(P/G).
$$
\item
If $P^\sharp$ is fine, every fractional ideal of $P$ is finitely
generated.
\end{enumerate}
\end{lemma}
\begin{proof} The first assertion is left to the reader. Next, we
remark that $I^{-1}/G=(I/G)^{-1}$ and $(IJ)/G=(I/G)\cdot(J/G)$,
for every fractional ideals $I,J$ of $P$, which imply immediately
assertions (ii) and (iii). Lastly, suppose that $P^\sharp$ is
finitely generated, and let $I$ be any fractional ideal of $P$;
pick $x\in I^{-1}$; since $P$ is integral, $I$ is finitely
generated if and only if the same holds for $xI$. Hence, in order
to show (iv), we may assume that $I\subset P$, in which case the
assertion follows from proposition \ref{prop_ideals-in-fg-mon}(ii)
and lemma \ref{lem_spec-quots}(i.a).
\end{proof}

In order to characterize the monoids $P$ such that $\Div(P)$
is a group, we make the following :

\begin{definition} Let $P$ be an integral monoid, and
$a\in P^\gp$ any element.
\begin{enumerate}
\alphaenu
\item
We say that $a$ is {\em power-bounded}, if there exists
$b\in P$ such that $a^nb\in P$ for all $n\in\N$.
\item
We say that $P$ is {\em completely saturated}, if all
power-bounded elements of $P^\gp$ lie in $P$.
\end{enumerate}
\end{definition}

\begin{example}\label{ex_ordered-gps}
Let $(\Gamma,\leq)$ be an ordered abelian group, and set
$\Gamma^+:=\{\gamma\in\Gamma~|~\gamma\leq 1\}$. Then
$\Gamma^+$ is always a saturated monoid, but it is
completely saturated if and only if the convex rank of
$\Gamma$ is $\leq 1$ (see \cite[Def.6.1.20]{Ga-Ra}). The
proof shall be left as an exercise for the reader.
\end{example}

\begin{proposition}\label{prop_completely-sat}
Let $P$ be an integral monoid. We have :
\begin{enumerate}
\item
$(\Div(P),\odot)$ is an abelian group if and only if $P$
is completely saturated.
\item
If $P$ is fine and saturated, then $P$ is completely saturated.
\item
Let $A$ be a Krull domain, and set $A':=A\!\setminus\!\{0\}$.
Then $(A',\cdot)$ is a completely saturated monoid.
\end{enumerate}
\end{proposition}
\begin{proof}(i): Suppose that $I\in\Div(I)$ admits an inverse
$J$ in the monoid $(\Div(P),\odot)$, and notice that
$I\odot I^{-1}\subset P$; it follows easily that
$I\odot(J\cup I^{-1})^*=P$, hence $I^{-1}\subset J$, by
the uniqueness of the inverse. On the other hand, if $J$
strictly contains $I^{-1}$, then $IJ$ strictly contains
$P$, which is absurd. Thus, we see that $\Div(P)$ is a
group if and only if $I\odot I^{-1}=P$ for every $I\in\Div(P)$.
Now, suppose first that $P$ is completely saturated. In
view of lemma \ref{lem_associat}(ii), we are reduced to
showing that $P$ is contained in every invertible fractional
ideal containing $I^{-1}I$. Hence, say that $I^{-1}I\subset aP$
for some $a\in P^\gp$; equivalently, we have
$a^{-1}I^{-1}I\subset P$, {\em i.e.} $a^{-1}I^{-1}\subset I^{-1}$,
and then $a^{-k}I^{-1}\subset I^{-1}$ for every integer $k\in\N$.
Say that $b\in I^{-1}$ and $c\in I^*$; we conclude that
$a^{-k}bc\in P$ for every $k\in\N$, so $a^{-1}\in P$, by
assumption, and finally $P\subset aP$, as required.

Conversely, suppose that $\Div(P)$ is a group, and let
$a\in P^\gp$ be any power-bounded element. By definition,
this means that the $P$-submodule $I$ of $P^\gp$ generated
by $(a^k~|~k\in\N)$ is a fractional ideal of $P$. Then
$I^{-1}$ is a reflexive fractional ideal, and by assumption
$I^{-1}$ admits an inverse, which must be $I^*$, by the
foregoing. On the other hand, by construction we have
$aI\subset I$, hence $aI^*=(aI)^*\subset I^*$. We deduce that
$aP=a(I^*\odot I^{-1})=aI^*\odot I^{-1}\subset I^*\odot I^{-1}=P$,
{\em i.e.} $a\in P$, as stated.

(ii) is a special case of the following :

\begin{claim}
Let $P$ be any fine and saturated monoid, $M\subset P^\gp$ a
non-empty finitely generated $P$-submodule, and $a\in P^\gp$
an element such that $aM\subset M$. Then $a\in P$.
\end{claim}
\begin{pfclaim} Pick any $m\in M$, and denote by $M'\subset M$
the submodule generated by $(a^km~|~k\in\N)$. According to
proposition \ref{prop_ideals-in-fg-mon}(i), there exists
$N\geq 0$ such that $M'$ is generated by the finite system
$(a^km~|~k=0,\dots,N)$. Especially, $a^{N+1}m\in M'$, and
therefore there exist $x\in P$ and $i\leq N$ such that
$a^{N+1}m=a^imx$ in $M$; it follows that $a^{N+1-i}\in P$,
and finally $a\in P$, since $P$ is saturated.
\end{pfclaim}

(iii): See \cite[\S12]{Mat} for the basic generalities
on Krull domains. One is immediately reduced to the case where
$A$ is a valuation ring whose value group $\Gamma$ has
rank $\leq 1$. Taking into account (i) and lemma
\ref{lem_rflx_and_quot}, it then suffices to show that the
monoid $A'/A^\times$ is completely sataurated. However,
the latter is isomorphic to the submonoid $\Gamma^+$ of
elements $\leq 1$ in $\Gamma$, so the assertion follows from
example \ref{ex_ordered-gps}.
\end{proof}

\sset\subsubsection{}\label{subsec_frac-ideal-mons}
Let $\phi:P\to Q$ be a morphism of integral monoids, and
$I$ any fractional ideal of $P$; notice that
$IQ:=\phi^\gp(I)Q\subset Q^\gp$ is a fractional ideal of $Q$.
Moreover, the identities
$$
(I_1\cup I_2)Q=I_1Q\cup I_2Q
\qquad
(I_1I_2)Q=(I_1Q)\cdot(I_2Q)
\qquad
\text{for all fractional ideals $I_1,I_2\subset P^\gp$}
$$
are immediate from the definitions. Likewise, if $A$ an
integral domain and $\alpha:P\to(A\!\setminus\!\{0\},\cdot)$
a morphism of monoids, then the $A$-submodule
$IA:=\alpha^\gp(I)A$ of the field of fractions of $A$ is a
fractional ideal of the ring $A$ (in the usual commutative
algebraic meaning : see example \ref{ex_integral-doms}),
and we have corresponding identities :
$$
(I_1\cup I_2)A=I_1A+I_2A
\qquad
(I_1I_2)A=(I_1A)\cdot(I_2A)
\qquad
\text{for all fractional ideals $I_1,I_2\subset P^\gp$}.
$$

\begin{lemma}\label{lem_rflx-rflx}
In the situation of \eqref{subsec_frac-ideal-mons}, suppose
that $\phi$ is flat and $A$ is $\alpha$-flat, and let
$I,J,J'\subset P^\gp$ be three fractional ideals, with $I$
finitely generated. Then we have :
\begin{enumerate}
\item
$(J:I)Q=(JQ:IQ)$ and $(J:I)A=(J:I)A$.
\item
Especially, if $I$ is reflexive, the same holds for $IQ$ and $IA$
(see example {\em\ref{ex_integral-doms}}).
\item
Suppose furthermore that $A$ is local, and $\alpha$ is a local
morphism. Then $JA=J'A$ if and only if $J=J'$.
\item
If $P$ is fine, the rule $I\mapsto IQ$ and $I\mapsto IA$ define
morphisms of monoids
$$
\Div(\phi):\Div(P)\to\Div(Q)
\qquad
\Div(\alpha):\Div(P)\to\Div(A)
$$
(where $\Div(A)$ is defined as in example {\em\ref{ex_integral-doms}}),
and $\Div(\alpha)$ is injective, if $\alpha$ is local and $A$ is a
local domain.
\end{enumerate}
\end{lemma}
\begin{proof}(i): Say that $I=a_1P\cup\cdots\cup a_nP$ for elements
$a_1,\dots,a_n\in P^\gp$. Then
$$
(J:I)=a_1^{-1}J\cap\cdots\cap a_n^{-1}J
\qquad\text{and}\qquad
(JQ:IQ)=a_1^{-1}JQ\cap\cdots\cap a_n^{-1}JQ
$$
and likewise for $(J:I)A$, hence the assertion follows from an easy
induction, and the following
\begin{claim}\label{cl_old-claim}
For any two fractional ideals $J_1,J_2\subset P$, we have
$(J_1\cap J_2)Q=J_1Q\cap J_2Q$ and $(J_1\cap J_2)A=J_1A\cap J_2A$.
\end{claim}
\begin{pfclaim} Pick any $x\in P$ such that $xJ_1,xJ_2\subset P$;
since $P$ is an integral monoid, and $A$ is an integral domain, it
suffices to show that $x(J_1\cap J_2)Q=xJ_1Q\cap xJ_2Q$ and likewise
for $x(J_1\cap J_2)A$, and notice that $x(J_1\cap J_2)=xJ_1\cap xJ_2$.
We may thus assume that $J_1$ and $J_2$ are ideals of $P$, in which
case the assertion is lemma \ref{lem_intersect-ideals}.
\end{pfclaim}

(ii): Suppose that $I$ is reflexive; from (i) we deduce that
$((IA)^{-1})^{-1}=IA$. The assertion is an immediate consequence,
once one remarks that, for any fractional ideal $J\subset A$, there
is a natural isomorphism of $A$-modules : $J^{-1}\isom
J^\vee:=\Hom_A(J,A)$. Indeed, the isomorphism assigns to any $x\in
J^{-1}$ the map $\mu_x:J\to A$ : $a\mapsto xa$ for every $a\in J$
(details left to the reader).

(iii): We may assume that $JA=J'A$, and we prove that $J=J'$,
and by replacing $J'$ by $J\cup J'$, we may assume that $J\subset J'$.
Then the contention follows easily from lemma \ref{lem_faithful-phi-flat}.

(iv): This is immediate from (i) and (iii).
\end{proof}

\begin{remark}\label{rem_pretty-obvious}
In the situation of lemma \ref{lem_rflx-rflx}(iv), obviously
$\Div(\phi)$ restricts to a morphism of submonoids :
$$
\Div_+(\phi):\Div_+(P)\to\Div_+(Q).
$$
\end{remark}

\sset\subsubsection{}
Next, suppose that $P$ is fine and saturated, and let
$I\subset P^\gp$ be any fractional ideal. Then theorem
\ref{th_structure-of-satu}(i) easily implies that :
$$
I^{-1}=\bigcap_{\hgt\,\fp=1}(I_\fp)^{-1}
$$
where the intersection -- running over the prime ideals of $P$ of
height one -- is taken within $\Hom_P(I,P^\gp)$, which naturally
contains all the $(I_\fp)^{-1}$. The structure of the fractional
ideals of $P_\fp$ when $\hgt\,\fp=1$ is very simple : quite
generally, theorem \ref{th_structure-of-satu}(ii) easily implies
that if $\dim P=1$, then all fractional ideals are cyclic, and
then clearly they are reflexive. On the other hand, $I$ is
finitely generated, by lemma \ref{lem_rflx_and_quot}(iv). We
deduce that $I$ is reflexive if and only if :
\set\begin{equation}\label{eq_reflex-on-fs-mon}
I=\bigcap_{\hgt\,\fp=1}I_\fp.
\end{equation}
Indeed, suppose that \eqref{eq_reflex-on-fs-mon} holds; then
we have $I^*=\bigcap_{\hgt\,\fp=1}(I_\fp^{-1})^{-1}=
\bigcap_{\hgt\,\fp=1}I_\fp$, since we have just seen that $I_\fp$
is a reflexive fractional ideal of $P_\fp$, for every
prime ideal $\fp$ of height one.

\begin{proposition}\label{prop_classify-reflex-mon}
Let $P$ be a fine and saturated monoid, and denote by
$D\subset\Spec\,P$ the subset of all prime ideals of height
one. Then the mapping :
\set\begin{equation}\label{eq_divisor-grp-mons}
\Z^{\oplus D}\to\Div(P)
\quad : \quad
\sum_{\hgt\,\fp=1}n_\fp[\fp]\mapsto
\bigcap_{\hgt\,\fp=1}\fm_{P_\fp}^{n_\fp}
\end{equation}
is an isomorphism of abelian groups.
\end{proposition}
\begin{proof} Here $\fm_{P_\fp}\subset P_\fp$ is the maximal ideal,
and for $n\geq 0$, the notation $\fm_{P_\fp}^n$ means the usual
$n$-th power operation in the monoid $\cP(P^\gp)$, which we extend
to all integers $n$, by letting $\fm_{P_\fp}^n:=\fm_{P_\fp}^{-n}$
whenever $n<0$.

In order to show that \eqref{eq_divisor-grp-mons} is well defined,
set $I:=\bigcap_{\hgt\,\fp=1}\fm_\fp^{n_\fp}$. Pick, for every $\fp$
such that $n_\fp<0$, an element $x_\fp\in\fp$, and set
$y_\fp:=x_\fp^{-n_\fp}$; if $n_\fp\geq 0$, set $y_\fp:=1$. Then it
is easy to check (using theorem \ref{th_structure-of-satu}(i)) that
$\prod_{\hgt\,\fp=1}y_\fp$ lies in $I^{-1}$, hence $I$ is a
fractional ideal. Next, for given $\fp,\fp'\in D$, notice that
$(P_\fp)_{\fp'}=P^\gp$; it follows that
\set\begin{equation}\label{eq_localiez-at-hgt-1}
I_\fp=\fm_\fp^{n_\fp} \qquad \text{for every $\fp\in D$}
\end{equation}
therefore $I$ is reflexive. Furthermore, it is easily seen (from
theorem \ref{th_structure-of-satu}(ii)), that every reflexive ideal
of $P_\fp$ is of the form $\fm_P^n$ for some integer $n$, and
moreover $\fm_P^n=\fm_P^m$ if and only if $n=m$. Then
\eqref{eq_reflex-on-fs-mon} implies that the mapping
\eqref{eq_divisor-grp-mons} is surjective, and the injectivity
follows from \eqref{eq_localiez-at-hgt-1}. It remains to check
that \eqref{eq_divisor-grp-mons} is a group homomorphism, and
to this aim we may assume -- in view of lemma
\ref{lem_rflx-rflx}(iv) -- that $\dim P=1$, in which case the
assertion is immediate.
\end{proof}

\sset\subsubsection{}\label{subsec_J_P-plus}
A morphism $\phi:I\to J$ of fractional ideals of $P$ is, by
definition, a morphism of $P$-modules. Let $x,y\in I$ be any
two elements; we may find $a,b\in P$ such that $ax=by$ in $I$,
and therefore $a\phi(x)=\phi(ax)=\phi(by)=b\phi(y)$; thus,
$\phi(y)=(b^{-1}a)\cdot\phi(x)=(x^{-1}y)\cdot\phi(x)$. This
shows that, for every morphism $\phi:I\to J$ of fractional
ideals, there exists $c\in P^\gp$ such that $\phi(x)=cx$ for
every $x\in I$.  Especially, $I\simeq J$ if and only if there
exists $a\in P^\gp$ such that $I=aJ$. Likewise one may
characterize the morphisms and isomorphisms of fractional
ideals of an integral domain. We denote
$$
\bar\Div(P)
$$
the set of isomorphism classes of reflexive fractional
ideals of $P$. From the foregoing, it is clear that if
$I\simeq I'$, we have $I\odot J\simeq I'\odot J$ for every
$J\in\Div(P)$; therefore the composition law of $\Div(P)$
descends to a composition law for $\bar\Div(P)$, which makes
it into a (commutative) monoid, and if $P$ is completely
saturated, then $\bar\Div(P)$ is an abelian group. We also
deduce an exact sequence of monoids
\set\begin{equation}\label{eq_basic-Div-seq}
1\to P^\times\to P^\gp\xrightarrow{\ j_P\ }\Div(P)\to
\bar\Div(P)\to 1
\end{equation}
where $j_P$ is given by the rule $A$: $a\mapsto aP$ for every
$a\in P$; especially, $j_P$ restricts to a morphism of monoids
$$
j^+_P:P\to\Div_+(P).
$$
Likewise, we define $\bar\Div(A)$, for any integral domain
$A$ : see example \ref{ex_integral-doms}. Moreover, in
the situation of \eqref{subsec_frac-ideal-mons}, we
see from lemma \ref{lem_rflx-rflx}(iv) that if $\alpha$
is local, $P$ is fine, $A$ is $\alpha$-flat, and $\phi$
is flat, then $\Div(\phi)$ and $\Div(\alpha)$ descend to
well defined morphisms of monoids
$$
\bar\Div(\phi):\bar\Div(P)\to\bar\Div(Q)
\qquad
\bar\Div(\alpha):\bar\Div(P)\to\bar\Div(A).
$$

\begin{proposition}\label{prop_I-I}
Let $P$ be a fine and saturated monoid, $I,J\subset P^\gp$
two fractional ideals, $A$ a local integral domain, and
$\alpha:P\to(A,\cdot)$ a local morphism of monoids. We have :
\begin{enumerate}
\item
$(I:I)=P$.
\item
Suppose that $A$ is $\alpha$-flat. Then $IA\simeq JA$ if and
only if $I\simeq J$. Especially, in this case $\bar\Div(\alpha)$
is an injective map.
\end{enumerate}
\end{proposition}
\begin{proof} (i): Clearly it suffices to show that
$(I:I)\subset P$. Hence, say that $x\in(I:I)$, and pick any
$a\in I$; it follows that $x^na\in P$ for every $n>0$; in
the additive group $\log P^\gp$ we have therefore the identity
$n\cdot\log(x)+\log(a)\in\log P$, so
$\log(x)+n^{-1}\log(a)\in(\log P)_\R$ for every $n>0$. Since
$(\log P)_\R$ is a convex polyhedral cone in $(\log P^\gp)_\R$,
we deduce that $x\in(\log P)_\R\cap(\log P^\gp)=\log P$
(proposition \ref{prop_Gordon}(iii)), as claimed.

(ii): We may assume that $IA$ is isomorphic to $JA$, and
we show that $I$ is isomorphic to $J$. Indeed, the assumption
means that $a(IA)=JA$ for some $x\in\Frac(A)$; therefore,
$a\in(JA:IA)$ and $a^{-1}\in(IA:JA)$, so
$$
A=(IA:JA)\cdot(JA:IA)=((I:J)\cdot(J:I))A
$$
by virtue of lemma \ref{lem_rflx-rflx}(i). Since $A$ is local,
it follows that there exist $a\in(I:J)$ and $b\in(J:I)$ such
that $\alpha(ab)\in A^\times$, whence $ab\in P^\times$, since
$\alpha$ is local. It follows easily that $I=aJ$, as asserted.
\end{proof}

\begin{example}\label{ex_explicit-dim-2}
(i)\ \
Let $P$ be a fine and saturated monoid, and $D\subset\Spec\,P$
the subset of all prime ideals of height one; for every
$\fp\in D$, denote
$$
v_\fp:P\to P_\fp^\sharp\isom\N
$$
the composition of the localization map, and the natural
isomorphism resulting from theorem \ref{th_structure-of-satu}(ii).
A simple inspection shows that the isomorphism
\eqref{eq_divisor-grp-mons} identifies the map $j_P$ of
\eqref{eq_basic-Div-seq} with the morphism of monoids
$$
v_P:P^\gp\to\Z^{\oplus D}
\qquad
x\mapsto(v^\gp_\fp(x)~|~\fp\in D).
$$
With this notation, the isomorphism \eqref{eq_divisor-grp-mons}
is the map given by the rule :
$$
k_\bullet\mapsto v_P^{-1}(k_\bullet+\N^{\oplus D})
\qquad
\text{for every $k_\bullet\in\Z^{\oplus D}$}.
$$

(ii)\ \
Suppose now that $P$ is sharp and $\dim P=2$, in which case
$D=\{\fp_1,\fp_2\}$ contains exactly two elements. According
to example \ref{ex_satu-dim-two}(ii), we may find a basis
$f_1,f_2$ of $P^\gp$, such that the two facets $P\setminus\fp_1$
and $P\setminus\fp_2$ of $P$ are generated respectively
by $e_1:=f_1$ and $e_2:=af_1+bf_2$, for some $a,b\in\N$,
with $a<b$ and $(a,b)=1$. It follows easily that $P$ is
a submonoid of the free monoid
$$
Q:=\N e'_1\oplus\N e'_2
\qquad
\text{where $e'_1:=b^{-1}e_1$ and $e'_2:=b^{-1}e_2$}
$$
and $Q^\gp/P^\gp\simeq\Z/b\Z$ (details left to the reader).
The induced map $\Spec\,Q\to\Spec\,P$ is a homeomorphism;
especially $Q$ admits two prime ideals $\fq_1$, $\fq_2$
of height one, so that $\fq_i\cap P=\fp_i$ for $i=1,2$,
whence -- by proposition \ref{prop_classify-reflex-mon} --
a natural isomorphism
$$
s^*:\Div(Q)\isom\Div(P)
$$
and notice that $j_Q:Q^\gp\to\Div(Q)$ is the isomorphism
given by the rule : $e'_i\mapsto\fq_i$ for $i=1,2$.
Moreover, we have commutative diagrams of monoids :
$$
{\diagram P \ar[r]^-s \ar[d]_{v_{\fp_i}} & Q \ar[d]^{v_{\fq_i}} \\
\N \ar[r]^-{t_i} & \N
\enddiagram}
\qquad
(i=1,2).
$$
Clearly, $Q\setminus\fq_i$ is the facet generated by $e'_i$,
so $v_{\fq_i}$ is none else than the projection onto the
direct factor $\N e'_{3-i}$, for $i=1,2$. In order to compute
$v_{\fp_i}$, it then suffices to determine $t_i$, or equivalently
$t^\gp_i$. However, set $\tau_i:=v^\gp_{\fq_i}\circ s^\gp$;
clearly $\tau_1(f_2)=\tau_1(e'_2-ae'_1)=1$, so $\tau_1$ is
surjective. Also, $\tau_2(f_1)=b$ and $\tau_2(f_2)=-a$, so
$\tau_2$ is surjective as well; therefore both $t_1$ and
$t_2$ are the identity endomorphism of $\N$. Summing up,
we find that
$$
j_P=s^*\circ j_Q\circ s^\gp
$$
and the morphism $v_P$ is naturally identified with
$s^\gp:P^\gp\to Q^\gp$. Especially, we have obtained a
natural isomorphism
$$
\bar\Div(P)\isom\Z/b\Z.
$$
We may then rephrase in more intrinsic terms the classification
of example \ref{ex_satu-dim-two}(iii) : namely the isomorphism
class of $P$ is completely determined by the datum of $\bar\Div(P)$
and the equivalence class of the height one prime ideals of
$P$ in the quotient set $\bar\Div(P)^\dagger$ defined as in
{\em loc.cit.}

(iii)\ \
In the situation of (ii), a simple inspection yields the
following explicit description of all reflexive fractional
ideals of $P$. Recall that such ideals are of the form
$$
I_{k_1,k_2}:=\fm_1^{k_1}\cap\fm_2^{k_2}=
\{x\in P^\gp~|~v_{\fp_1}(x)\geq k_1,\ v_{\fp_2}(x)\geq k_2\}
$$
where $\fm_i$ is the maximal ideal of $P_{\fp_i}$, and
$k_i\in\Z$, for $i=1,2$. Then
$$
I_{k_1,k_2}=\{x_1e_1+x_2e_2~|~
x_1,x_2\in b^{-1}\Z,\ x_1\geq b^{-1}k_2,\ x_2\geq b^{-1}k_1\}
\cap P^\gp
\qquad
\text{for all $k_1,k_2\in\Z$}.
$$
With this notation, the cyclic reflexive ideals are then those
of the form
$$
(x_1f_1+x_2f_2)P=I_{x_2b,x_1b-x_2a}
\qquad
\text{with $x_1,x_2\in\Z$}.
$$
Especially, we see that the classes of $\fp_1=I_{0,1}$ and
$\fp_2=I_{1,0}$ are both of order $b$ in $\bar\Div(P)$.
\end{example}

The following estimate, special to the two-dimensional case,
will be applied to the proof of the almost purity theorem for
towers of regular log schemes.

\begin{lemma} Let $P$ be a fine and saturated monoid of dimension
$2$, and denote by $b$ the order of the finite cyclic group
$\bar\Div(P)$. We have :
$$
\fm_P^{[b/2]}\subset I\cdot I^{-1}
\qquad
\text{for every $I\in\Div(P)$}
$$
(where $[b/2]$ denotes the largest integer $\leq b/2$).
\end{lemma}
\begin{proof} Notice first that the assertion holds for a
given $I\in\Div(P)$, if and only if it holds for $xI$, for
any $x\in P^\gp$. If $b=1$, then $P=\N^{\oplus 2}$, in
which case $\bar\Div(P)=0$, so every reflexive fractional ideal
of $P$ is isomorphic to $P$, and the assertion is clear.
Hence, assume that $b>1$; let $\fp_1$, $\fp_2$ be the two
prime ideals of height one of $P$, and define $Q$, $e_1$,
$e_2$ and $I_{k_1,k_2}$ for every $k_1,k_2\in\Z$, as in
example \ref{ex_explicit-dim-2}(ii,iii). With this notation,
notice that
$$
\fm_P\setminus\fm_P^2=\{e_1,e_2\}\cup\Sigma
\qquad
\text{where $\Sigma\subset
\{x_1e_1+x_2e_2~|~x_1,x_2\in b^{-1}\Z,\ 0\leq x_1,x_2<1\}$}.
$$
It follows easily that, for every $i\in\N$, every element of
$\fm^i_P$ is of the form $x_1e_1+x_2e_2$ with
$x_1,x_2\in b^{-1}\N$ and $\max(x_1,x_2)\geq b^{-1}i$.
Hence, let $I\in\Div(P)$ and $x:=x_1e_1+x_2e_2\in\fm_P^{[b/2]}$,
and say that $bx_1\geq[b/2]$. According to example
\ref{ex_explicit-dim-2}(iii), we may assume that
$I=\fp_2^j=I_{j,0}$ for some $j\in\{0,\dots,b-1\}$.
Moreover, notice that the assertion holds for $I$ if and
only if it holds for $I^{-1}$, whose class in $\bar\Div(P)$
agrees with the class of $\fp_2^{b-j}$. Clearly, either
$j\leq[b/2]$ or $b-j\leq[b/2]$; hence, we may assume that
$j\in\{0,\dots,[b/2]\}$. Thus, $P\subset I^{-1}$, and
$I\subset I\cdot I^{-1}$, and clearly $x\in I$, so we are done
in this case. The case where $bx_2\geq[b/2]$ is dealt with in
the same way, by writing $I=\fp_1^j$ for some non-negative
$j\leq[b/2]$ : the details are left to the reader.
\end{proof}

If $f:P\to Q$ is a general morphism of integral monoids, and
$I$ a fractional ideal of $Q$, the $P$-module $f^{\gp-1}(I)$
is not necessarily a fractional ideal of $P$ (for instance,
consider the natural map $P\to P^\gp$). One may obtain some
positive results, by restricting to the class of morphisms
introduced by the following :

\begin{definition}\label{def_Kummer-monoids}
Let $f:P\to Q$ be a morphism of monoids. We say that $f$ is
{\em of Kummer type}, if $f$ is injective, and the induced
map $f_\Q:P_\Q\to Q_\Q$ is surjective (notation of
\eqref{subsec_from-con-to-mon}).
\end{definition}

\begin{lemma}\label{lem_Kummer-fans}
Let $f:P\to Q$ be a morphism of monoids of Kummer type,
$S_Q\subset Q$ a submonoid, and set $S_P:=f^{-1}S_Q$. We have :
\begin{enumerate}
\item
The map\/ $\Spec\,f:\Spec\,Q\to\Spec\,P$ is bijective; especially
$\dim P=\dim Q$.
\item
If $Q^\times$ is a torsion-free abelian group, $P$ is the trivial
monoid (resp. is sharp) if and only if the same holds for $Q$.
\item
The induced morphism $S_P^{-1}P\to S_Q^{-1}Q$ is of Kummer type.
\item
If $P$ is integral, the unit of adjunction $P\to P^\sat$ is of
Kummer type.
\item
Suppose that $P$ is integral and saturated. Then
$f^\sharp:P^\sharp\to Q^\sharp$ is of Kummer type.
\item
If both $P$ and $Q$ are integral, and $P$ is saturated, then
$f$ is exact.
\end{enumerate}
\end{lemma}
\begin{proof} (ii) and (iv) are trivial, and (iii) is an exercise for the
reader.

(i): Let $F,F'\subset Q$ be two faces such that $f^{-1}F=f^{-1}F'$,
and say that $x\in F$. Then $x^n\in f(P)$ for some $n>0$, so $x^n\in
f(f^{-1}F')$, whence $x\in F'$, which implies that $\Spec\,f$ is
injective. Next, for a given face $F$ of $P$, let $F'\subset Q$ be
the subset of all $x\in Q$ such that there exists $n>0$ with $x^n\in
f(F)$. It is easily seen that $F'$ is a face of $Q$, and moreover
$f^{-1}F'=F$, which shows that $\Spec\,f$ is also surjective.

(v): Clearly the map $(P^\sharp)_\Q\to(Q^\sharp)_\Q$ is
surjective. Now, let $x,y\in P$ such that the images of $f(x)$
and $f(y)$ agree in $Q^\sharp$, {\em i.e.} there exists
$u\in Q^\times$ with $u\cdot f(x)=f(y)$; we may find $n>0$ such
that $u^n,u^{-n}\in f(P)$. Say that $u^n=f(v)$, $u^{-n}=f(w)$;
since $f(vw)=1$, we have $vw=1$, and moreover $f(vx^n)=f(y^n)$,
so $vx^n=y^n$. Therefore $x^ny^{-n},x^{-n}y^n\in P$, and since
$P$ is saturated we deduce that $xy^{-1},x^{-1}y\in P$, so the
images of $x$ and $y$ agree in $P^\sharp$.

(vi): Notice first that $f^\gp$ is injective, since the same
holds for $f$. Suppose $x\in P^\gp$ and $f^\gp(x)\in Q$; we
may then find an integer $k>0$ and $y\in P$ such that
$f(y)=f(x)^k$. Since $P$ is saturated, it follows that
$x\in P$, so $f$ is exact.
\end{proof}

\sset\subsubsection{}\label{subsec_one-sat-one-not}
Suppose that $\phi:P\to Q$ is a morphism of integral monoids
of Kummer type, with $P$ saturated, and let $I\subset Q^\gp$
be a fractional ideal.
Then $\phi^*I:=\phi^{\gp-1}(I)$ is a fractional ideal of $P$.
Indeed, by assumption there exists $a\in Q$ such that
$aI\subset Q$; we may find $k>0$ and $b\in P$ such that
$a^k=\phi(b)$, so $\phi(bx)\in\phi(P)^\gp\cap Q=\phi(P)$ for
every $x\in\phi^*I$, since $\phi$ is exact (lemma
\ref{lem_Kummer-fans}(v)); therefore $b\cdot\phi^*(I)\subset P$.

\sset\subsubsection{}\label{subsec_ramif-index}
In the situation of \eqref{subsec_one-sat-one-not}, suppose
that both $P$ and $Q$ are fine and saturated, and let
$\gr_\bullet Q^\gp$ be the $\phi$-grading of $Q$, indexed by
$(\Gamma,+):=Q^\gp/P^\gp$ (see remark \ref{rem_why-not}(iii));
for every $x\in Q^\gp$, denote $\bar x\in\Gamma$ the image of
$x$. Let also $I\subset Q^\gp$ be any fractional ideal, and
denote by $\gr_\bullet I$ the $\Gamma$-grading on $I$ deduced
from the $\phi$-grading of $Q^\gp$; arguing as in
\eqref{subsec_one-sat-one-not}, it is easily seen that, more
generally, $\phi^*(x^{-1}\gr_{\bar x}I)$ is a fractional
ideal of $P$, for every $x\in Q$ (details left to the reader).
For every prime ideal $\fq$ of height one in $Q$, we have
a commutative diagram of monoids :
\set\begin{equation}\label{eq_ramif-index}
{\diagram
P \ar[r]^-\phi \ar[d]_{v_{\phi^{-1}\fq}} & Q \ar[d]^{v_\fq} \\
\N \ar[r]^{e_{\fq}} & \N
\enddiagram}
\end{equation}
where $v_\fq$ and $v_{\phi^{-1}\fq}$ are defined as in example
\ref{ex_explicit-dim-2}(i), and $e_\fq$ is the multiplication
by a non-zero (positive) integer, which we call the
{\em ramification index of\/ $\phi$ at\/ $\fq$}, and we denote
also $e_\fq$.

\begin{lemma}\label{lem_grad-reflex}
In the situation of \eqref{subsec_ramif-index}, suppose that
$I$ is a reflexive fractional ideal. Then
$\phi^*(a^{-1}\cdot\gr_{\bar a}I)$ is a reflexive fractional
ideal of $P$, for every $a\in Q^\gp$.
\end{lemma}
\begin{proof} Clearly, we may replace $I$ by $a^{-1}I$, and
reduce to the case where $a=1$, in which case we have to check
that $\phi^*I$ is a reflexive fractional ideal. However,
according to example \ref{ex_explicit-dim-2}(i), we may write
$I=v_Q^{-1}(k_\bullet+\N^{\oplus D})$, where $D\subset\Spec\,Q$
is the subset of the height one prime ideals, and
$k_\bullet\in\Z^{\oplus D}$. Set
$$
k'_\bullet:=([e_\fq^{-1}k_\fq]~|~\fq\in D)
$$
(where, for a real number $x$, we let $[x]$ be the smallest
integer $\geq x$). Since $\Spec\,\phi$ is bijective (lemma
\ref{lem_Kummer-fans}(i)), the commutative diagrams
\eqref{eq_ramif-index} imply that
$\phi^*I=v_P^{-1}(k'_\bullet+\N^{\oplus D})$, whence the
contention.
\end{proof}

\begin{example}\label{ex_Kummer-dim-two}
Let $P$ be as in example \ref{ex_explicit-dim-2}(ii), set
$Q:=\Div_+(P)$, and take $\phi:=j^+_P:P\to Q$ (notation of
\eqref{subsec_J_P-plus}).
The discussion of {\em loc.cit.} shows that $\phi$ is a
morphism of Kummer type, and notice that the $\phi$-grading
of $Q$ is indexed by $Q^\gp/P^\gp=\bar\Div(P)$. Now, pick any
$x\in Q$, and let $\bar x\in\bar\Div(P)$ be the equivalence class
of $x$; by lemma \ref{lem_grad-reflex}, the $P$-module
$\gr_{\bar x}Q$ is isomorphic to a reflexive fractional ideal
of $P$. We claim that the isomorphism class of
$\gr_{\bar x}Q$ is precisely $\bar x{}^{-1}$ (where the inverse
is formed in the commutative group $\bar\Div(P)$). Indeed, let
$a\in P^\gp$ be any element; by definition, we have
$a\in\phi^*(x^{-1}\gr_{\bar x}Q)$ if and only if
$\phi^\gp(a)\in x^{-1}\gr_{\bar x}Q$, if and only if
$ax\in Q$, if and only if $a\in x^{-1}$, whence the claim.
Thus, the family
$$
(\gr_\gamma\Div_+(P)~|~\gamma\in\bar\Div(P))
$$
is a complete system of representatives for the isomorphism
classes of the reflexive fractional ideals of a fine, sharp
and saturated monoid $P$ of dimension $2$.
\end{example}

\begin{remark} Further results on reflexive fractional ideals
for monoids, and their divisor class groups can be found in
\cite{Ch}.
\end{remark}

\subsection{Fans}
According to Kato (\cite[\S9]{Ka2}), a fan is to a monoid what a
scheme is to a ring. More prosaically, the theory of fans is a
reformulation of the older theory of rational polyhedral
decompositions, developed in \cite{KKMS}.

\begin{definition}\label{def_mon-spaces}
(i)\ \ A {\em monoidal space} is a datum $(T,\cO_T)$ consisting
of a topological space $T$ and a sheaf of monoids $\cO_T$ on $T$.
\begin{enumerate}
\addenu
\item
A {\em morphism of monoidal spaces\/} is a datum
$$
(f,\log f):(T,\cO_T)\to(S,\cO_S)
$$
consisting of a continuous map $f:T\to S$, and a morphism $\log
f:f^*\cO_S\to\cO_T$ of sheaves of monoids that is {\em local}, {\em
i.e.} whose stalk $(\log f)_t:\cO_{S,f(t)}\to\cO_{T,t}$ is a local
morphism, for every $t\in T$. The {\em strict locus\/} of $(f,\log f)$
is the subset
$$
\mathrm{Str}(f,\log f)\subset T
$$
consisting of all $t\in T$ such that $(\log f)_t$ is an isomorphism.
\item
We say that a monoidal space $(T,\cO_T)$ is {\em sharp}, (resp. {\em
integral}, resp. {\em saturated}) if $\cO_T$ is a sheaf of sharp
(resp. integral, resp. integral and saturated) monoids.
\item
For any monoidal space (resp. integral monoidal space) $(T,\cO_T)$,
the {\em sharpening\/} (resp. the {\em saturation}) of $(T,\cO_T)$
is the sharp monoidal space $(T,\cO_T)^\sharp:=(T,\cO_T^\sharp)$
(resp. $(T,\cO_T)^\sat:=(T,\cO_T^\sat)$).
\end{enumerate}
\end{definition}

It is easily seen that the rule $(T,\cO_T)\mapsto(T,\cO_T)^\sharp$
extends to a functor from the category of monoidal spaces to the
full subcategory of sharp monoidal spaces. This functor is right
adjoint to the corresponding fully faithful embedding of categories.

Likewise, the functor $(T,\cO_T)\mapsto(T,\cO_T)^\sat$ is right
adjoint to the fully faithful embedding of the category of saturated
monoidal spaces, into the category of integral monoidal spaces.

\sset\subsubsection{}
Let $P$ be any monoid; for every $f\in P$, let us set
$$
D(f):=\{\fp\in\Spec\,P~|~f\notin\fp\}.
$$
Notice that $D(f)\cap D(g)=D(fg)$ for every $f,g\in P$. We endow
$\Spec\,P$ with the topology having a basis consisting of the
subsets $D(f)$, for every $f\in P$. Notice that $\fm_P$ is the
only closed point of $\Spec\,P$ (especially, $\Spec\,P$ is trivially
quasi-compact).

By lemma \ref{lem_std-verif}, the localization map $j_f:P\to P_f$
induces an identification $j^*_f:\Spec\,P_f\isom D(f)$. It is easily
seen that $(j^*_f)^{-1}D(fg)=D(j(g))\subset\Spec\,P_f$; in other
words, the topology of $\Spec\,P_f$ agrees with the topology induced
from $\Spec\,P$, via $j^*$.

Next, for every $f\in P$ we set :
$$
\cO_{\Spec\,P}(D(f)):=P_f.
$$
We claim that $\cO_{\Spec\,P}(D(f))$ depends only on the open subset
$D(f)$, up to natural isomorphism (and not on the choice of $f$).
More precisely, say that $D(f)\subset D(g)$ for two given elements
$f,g\in P$; it follows that the image of $g$ in $P_f$ lies outside
the maximal ideal $\fm_{P_f}$, hence $g\in P^\times_f$, and
therefore the localization map $j_f:P\to P_f$ factors uniquely
through a morphism of monoids :
$$
j_{f,g}:P_g\to P_f.
$$
Likewise, if $D(g)\subset D(f)$ as well, the localization $j_g:P\to
P_g$ factors through a unique map $j_{g,f}:P_f\to P_g$, whence the
identities :
$$
j_{f,g}\circ j_{g,f}\circ j_f=j_f \qquad j_{g,f}\circ j_{f,g}\circ
j_g=j_g
$$
and since $j_f$ and $j_g$ are epimorphisms, we see that $j_{f,g}$
and $j_{g,f}$ are mutually inverse isomorphisms.

\sset\subsubsection{}\label{subsec_def-Spec-P}
Say that $D(f)\subset D(g)\subset D(h)$ for some $f,g,h\in P$; by
direct inspection, it is clear that $j_{f,g}\circ j_{g,h}=j_{f,h}$,
so the rule $D(f)\mapsto P_f$ yields a well defined presheaf of
monoids on the site $\cC_P$ of open subsets of $\Spec\,P$ of the
form $D(f)$ for some $f\in P$. Then $\cO_{\Spec\,P}$ is trivially
a sheaf on $\cC_P$ (notice that if $D(f)=\bigcup_{i\in I}D(g_i)$
is an open covering of $D(f)$, then $D(g_i)=D(f)$ for some $i\in I$).
According to \cite[Ch.0, \S3.2.2]{EGAI} it follows that
$\cO_{\Spec\,P}$ extends uniquely to a well defined sheaf of monoids
on $\Spec\,P$, whence a monoidal space $(\Spec\,P,\cO_{\Spec\,P})$.
By inspecting the construction, we find natural identifications :
\set\begin{equation}\label{eq_ident-stalks}
(\cO_{\Spec\,P})_\fp\isom P_\fp \qquad \text{for every
$\fp\in\Spec\,P$}
\end{equation}
and moreover :
$$
P\isom\Gamma(\Spec\,P,\cO_{\Spec\,P}).
$$
It is also clear that the rule
\set\begin{equation}\label{eq_funct-with-no-name}
P\mapsto(\Spec\,P,\cO_{\Spec\,P})
\end{equation}
defines a functor from the category $\Mnd^o$ to the category of
monoidal spaces.

\begin{proposition}\label{prop_Spec-fan}
The functor \eqref{eq_funct-with-no-name} is right adjoint to the
functor :
$$
(T,\cO_T)\mapsto\Gamma(T,\cO_T)
$$
from the category of monoidal spaces, to the category $\Mnd^o$.
\end{proposition}
\begin{proof} Let $f:P\to\Gamma(T,\cO_T)$ be a map of monoids.
We define a morphism
$$
\phi_f:=(\phi_f,\log\phi_f):(T,\cO_T)\to(\Spec\,P,\cO_{\Spec\,P})
$$
as follows. Given $t\in T$, let $f_t:P\to\cO_{T,t}$ be the morphism
deduced from $f$, and denote by $\fm_t\subset\cO_{T,t}$ the maximal
ideal. We set $\phi_f(t):=f_t^{-1}(\fm_t)$. In order to show that
$\phi_f$ is continuous, it suffices to prove that
$U_s:=\phi_f^{-1}(D(s))$ is open in $M_T$, for every $s\in P$.
However, $U_s=\{t\in T~|~f_t(s)\in\cO_{T,t}^\times\}$, and it is
easily seen that this condition defines an open subset (details left
to the reader). Next, we define $\log\phi_f$ on the basic open
subsets $D(s)$. Indeed, let $j_s:\Gamma(T,\cO_T)\to\cO_T(U_s)$ be
the natural map; by construction, $j_s\circ f(s)$ is invertible in
$\cO_T(U_s)$, hence $j_s\circ f$ extends to a unique map of monoids :
$$
P_s=\cO_{\Spec\,P}(D(s))\to\phi_{f*}\cO_T(D(s)).
$$
By \cite[Ch.0, \S3.2.5]{EGAI}, the above rule extends to a unique
morphism $\cO_{\Spec\,P}\to\phi_{f*}\cO_T$ of sheaves of monoids,
whence -- by adjunction -- a well defined morphism
$\log\phi_f:\phi_f^*\cO_{\Spec\,P}\to\cO_T$. In order to show that
$(\phi_f,\log\phi_f)$ is the sought morphism of monoidal spaces, it
remains to check that $(\log\phi_f)_t:P_{\phi_f(t)}\to\cO_{T,t}$ is
a local morphism, for every $t\in T$. However, let $i_t:P\to
P_{\phi_f(t)}$ be the localization map; by construction, we have
$(\log\phi_f)_t\circ i_t=f_t$, and the contention is a
straightforward consequence.

Conversely, say that
$(\phi,\log\phi):(T,\cO_T)\to(\Spec\,P,\cO_{\Spec\,P})$ is a
morphism of monoidal spaces; then $\log\phi$ corresponds to a unique
morphism $\psi:\cO_{\Spec\,P}\to\phi_*\cO_T$, and we set
$$
f_\phi:=\Gamma(\Spec\,P,\psi):P\to\Gamma(T,\cO_T).
$$
By inspecting the definitions, it is easily seen that $f_{\phi_f}=f$
for every morphism of monoids $f$ as above. To conclude, it remains
only to show that the rule $(\phi,\log\phi)\mapsto f_\phi$ is
injective. However, for a given morphism of monoidal spaces
$(\phi,\log\phi)$ as above, and every point $t\in T$, we have a
commutative diagram of monoids :
$$
\xymatrix{ P \ar[r]^-{f_\phi} \ar[d] & \Gamma(T,\cO_T) \ar[d]^{j_t} \\
P_{\phi(t)} \ar[r]^-{\log\phi_t} & \cO_{T,t}}.
$$
Since $\log\phi_t$ is local, it follows that $\phi(t)=(f_\phi\circ
j_t)^{-1}\fm_t$, especially $f_\phi$ determines $\phi:T\to\Spec\,P$.
Finally, since the map $P\to P_{\phi(t)}$ is an epimorphism, we
see that $\log\phi_t$ is determined by $f_\phi$ as well, and the
proposition follows.
\end{proof}

\begin{definition} Let $(T,\cO_T)$ be a sharp monoidal space.
\begin{enumerate}
\item
We say that $(T,\cO_T)$ is an {\em affine fan}, if there exist a
monoid $P$ and an isomorphism of sharp monoidal spaces
$(\Spec\,P,\cO_{\Spec\,P})^\sharp\isom(T,\cO_T)$.
\item
In the situation of (i), if $P$ can be chosen to be finitely
generated (resp. fine), we say that $(T,\cO_T)$ is a {\em finite}
(resp. {\em fine}) {\em affine fan}.
\item
We say that $(T,\cO_T)$ is a {\em fan}, if there exists an open
covering $T=\bigcup_{i\in I}U_i$, such that the induced sharp
monoidal space $(U_i,\cO_{T|U_i})$ is an affine fan, for every $i\in
I$. We denote by $\Fan$ the full subcategory of the category of
monoidal spaces, whose objects are the fans.
\item
In the situation of (iii), if the covering $(U_i~|~i\in I)$ can be
chosen, so that $(U_i,\cO_{T|U_i})$ is a finite (resp. fine) affine
fan for every $i\in I$, we say that $(T,\cO_T)$ is {\em locally
finite} (resp. {\em locally fine}).
\item
We say that the fan $(T,\cO_T)$ is {\em finite} (resp. {\em fine})
if it is locally finite (resp. locally fine) and quasi-compact.
\item
Let $(T,\cO_T)$ be a fan. The {\em simplicial locus\/}
$T_\mathrm{sim}\subset T$ is the subset of all $t\in T$ such
that $\cO_{T,t}$ is a free monoid of finite rank.
\end{enumerate}
\end{definition}

\begin{remark}\label{rem_product-of-aff-fans}
(i)\ \  
For every monoid $P$, let $T_P$ denote the affine fan
$(\Spec\,P)^\sharp$. In light of proposition \ref{prop_Spec-fan},
it is easily seen that the functor $P\mapsto T_P$ is an equivalence
from the opposite of the full subcategory of sharp monoids, to the
category of affine fans.

(ii)\ \ Since the saturation functor commutes with localizations
(lemma \ref{lem_exc-satura}(i)), it is easily seen that the
saturation of a fan is a fan, and more precisely, the saturation
of an affine fan $T_P$, is naturally isomorphic to $T_{P^\sat}$.

(iii)\ \ Let $Q_1$ and $Q_2$ be two monoids; since the product
$P\times Q$ is also the coproduct of $P$ and $Q$ in the category
$\Mnd$ (see example \ref{ex_forrinst}(i)), we have a natural
isomorphism in the category of fans :
$$
T_{P\times Q}\isom T_P\times T_Q.
$$
More generally, suppose that $P\to Q_i$, for $i=1,2$, are two
morphisms of monoids. Then we have a natural isomorphism of fans :
$$
T_{Q_1\otimes_PQ_2}\isom T_{Q_1}\times_{T_P}T_{Q_2}.
$$
From this, a standard argument shows that fibre products are
representable in the category of fans. 

(iv)\ \ Furthermore, lemma \ref{lem_spec-quots}(ii) implies
that the natural map :
$$
\pi:\Spec\,(P\times Q)\to\Spec\,P\times\Spec\,Q
$$
is a homeomorphism (where the product of $\Spec\,P$ and $\Spec\,Q$
is taken in the category of topological spaces and continuous maps).

(v)\ \ Moreover, we have natural isomorphisms of monoids :
$$
\cO_{T_{P\times Q},{\pi^{-1}(s,t)}}\isom\cO_{T_P,s}\times\cO_{T_Q,t}
\qquad \text{for every $s\in\Spec\,P$ and $t\in\Spec\,Q$.}
$$

(vi)\ \ For any fan $T:=(T,\cO_T)$, and any monoid $M$, we shall
use the standard notation :
$$
T(M):=\Hom_\Fan((\Spec\,M)^\sharp,T).
$$
Especially, if $T$ is an affine fan, say $T=(\Spec\,P)^\sharp$,
then $T(M)=\Hom_\Mnd(P,M^\sharp)$; for instance, if $T$ is affine,
$T(\N)$ is a monoid, and $T(\Q_+)^\gp$ is a $\Q$-vector space.
Furthermore, by standard general nonsense we have natural
identifications of sets :
$$
(T_1\times_TT_2)(M)\isom T_1(M)\times_{T(M)}T_2(M)
$$
for any pair of $T$-fans $T_1$ and $T_2$, and every monoid $M$.
If $T$, $T_1$ and $T_2$ are affine, this identification is also
an isomorphism of monoids.
\end{remark}

\begin{example}\label{ex_fans}
(i)\ \  The topological space underlying the affine fan
$(\Spec\,\N,\cO_{\Spec\,\N})^\sharp$ consists of two points :
$\Spec\,\N=\{\emptyset,\fm\}$, where $\fm:=\N\!\setminus\!\{0\}$ is
the closed point. The structure sheaf $\cO:=\cO_{\Spec\,\N}$ is
determined as follows. The two stalks are $\cO_\emptyset=\{1\}$ (the
trivial monoid) and $\cO_\fm=\N$; the global sections are
$\Gamma(\Spec\,\N,\cO)=\N$.

(ii)\ \ Let $(T,\cO_T)$ be any fan, $P$ any monoid, with maximal
ideal $\fm_P$, and
$\phi:T_P:=(\Spec\,P,\cO_{\Spec\,P})^\sharp\to(T,\cO_T)$
a morphism of fans. Say that $\phi(\fm_P)\in U$ for some
affine open subset $U\subset F$; then
$\phi(\Spec\,P)\subset U$, hence $\phi$ factors through a
morphism of fans $T_P\to(U,\cO_{T|U})$. In view of proposition
\ref{prop_Spec-fan}, such a morphism corresponds to a unique
morphism of monoids $\phi^\natural:\cO_T(U)\to P^\sharp$,
and then
$\phi(\fm_P)=\phi^{\natural -1}(\fm_P)\in\Spec\,\cO_T(U)=U$.
The map on stalks determined by $\phi$ is the local morphism
$$
\cO_{T,\phi(\fm_P)}\isom
\cO_T(U)_{\phi(\fm_P)}/\cO_T(U)_{\phi(\fm_P)}^\times\to P^\sharp
$$
obtained from $\phi^\natural$ after localization at the prime ideal
$\phi^{\natural -1}(\fm_P)$.

(iii)\ \ For any two monoids $M$ and $N$, denote by
$\mathrm{loc.Hom}_\Mnd(M,N)$ the set of local morphisms of monoids
$M\to N$. The discussion in (ii) leads to a natural identification :
$$
T(P)\isom\coprod_{t\in T}\mathrm{loc.Hom}_\Mnd(\cO_{T,t},P^\sharp)
$$
For any monoid $P$. The {\em support\/} of a $P$-point $\phi\in T(P)$
is the unique point $t\in T$ such that $\phi$ corresponds to a
local morphism $\cO_{T,t}\to P^\sharp$.
\end{example}

\begin{example}\label{ex_multiply-by-n-in-fan}
(i)\ \ 
Let $P$ be any monoid, $k>0$ any integer, set
$T_P:=(\Spec\,P)^\sharp$, and let $\bek_P:P\to P$ be the
$k$-Frobenius map of $P$ (definition \ref{def_exact-phi}(ii)).
Then
$$
\bek_{T_P}:=\Spec\,\bek_P:T_P\to T_P
$$
is a well defined endomorphism inducing the identity on the
underlying topological space.

(ii)\ \ 
More generally, let $F$ be any fan; for every integer $k>0$
we have the {\em $k$-Frobenius endomorphism}
$$
\bek_F:F\to F
$$
which induces the identity on the underlying topological space, and
whose restriction to any affine open subfan $U\subset F$ is the
endomorphism $\bek_U$ defined as in (i).
\end{example}

\sset\subsubsection{}
Let $P$ be a monoid, $M$ a $P$-module, and set $T_P:=(\Spec\,P)^\sharp$.
We define a presheaf $M^\sim$ on the site of basic affine open subsets
$D(f)\subset\Spec\,P$ (for all $f\in P$), by the rule :
$$
U\mapsto M^\sim(U):=M\otimes_P\cO_{T_P}(U)
$$
(and for an inclusion $U'\subset U$ of basic open subsets, the
corresponding morphism $M^\sim(U)\to M^\sim(U')$ is deduced from
the restriction map $\cO_{T_P}(U)\to\cO_{T_P}(U')$).
It is easily seen that $M^\sim$ is a sheaf, hence it extends to
a well defined sheaf of $\cO_{T_P}$-modules on $T_P$
(\cite[Ch.0, \S3.2.5]{EGAI}). Clearly $\Gamma(T_P,M^\sim)=M$,
and the rule $M\mapsto M^\sim$ yields a well defined functor
$P\Mod\to\cO_{T_P}\Mod$, which is left adjoint to the global
section functor on $\cO_{T_P}$-modules :
$\cM\mapsto\Gamma(T_P,\cM)$ (verification left to the reader).

\begin{definition}\label{def_quasicoh-fans}
Let $(T,\cO_T)$ be a fan, $\cM$ a $\cO_T$-module.
We say that $\cM$ is {\em quasi-coherent}, if there exist an open
covering $T=\bigcup_{i\in I}U_i$ of $T$ by affine open subsets, and
for each $i\in I$ an $\cO_T(U_i)$-module $M_i$ with an isomorphism
of $\cO_{T|U_i}$-modules $\cM_{|U_i}\isom M^\sim_i$.
\end{definition}

\begin{remark}\label{rem_local-trivial-qcoh}
(i)\ \  
Let $(T,\cO_T)$ be a fan, $\cM$ a quasi-coherent $\cO_T$-module,
and $U\subset T$ be any open subset, such that $(U,\cO_{T|U})$
is an affine fan (briefly : an affine open subfan of $T$). Then,
since $U$ admits a unique closed point $t\in U$, it is easily
seen that $\cM_{|U}$ is naturally isomorphic to $\cM_t^\sim$ as
a $\cO_{T|U}$-module.

(ii)\ \ 
In the same vein, if $\cM$ is an invertible $\cO_T$-module (see
definition \ref{def_coh-idea-log}(iv)), then the restriction
$\cM_{|U}$ of $\cM$ to any affine open subset, is isomorphic to
$\cO_{T|U}$.

(iii)\ \
For any fan $(T,\cO_T)$, the sheaf of abelian groups $\cO^\gp_T$
is quasi-coherent (exercise for the reader). Suppose that $T$ is
integral; then an $\cO_T$-submodule $\cI\subset\cO^\gp_T$
is called a {\em fractional ideal\/} (resp. a {\em reflexive
fractional ideal}) of $\cO_T$ if $\cI$ is quasi-coherent, and
$\cI(U)$ is a fractional ideal (resp. a reflexive fractional
ideal) of $\cO_T(U)$, for every affine open subset $U\subset T$.

(iv)\ \
Let $P$ be any integral monoid, $I\subset P^\gp$ a fractional
ideal of $P$, and set $T_P:=(\Spec\,P)^\sharp$. It follows easily
from lemma \ref{lem_rflx_and_quot}(i) that
$I^\sim\subset\cO_{T_P}^\gp$ is a fractional ideal of $\cO_{T_P}$,
and $I^\sim$ is reflexive if and only if $I$ is a reflexive
fractional ideal of $P$ (lemma \ref{lem_rflx_and_quot}(ii)).

(v)\ \
Suppose that $T$ is locally finite; in this case, it follows
easily from proposition \ref{prop_ideals-in-fg-mon}(ii) that
every quasi-coherent ideal of $\cO_T$ is coherent. Likewise,
if $T$ is also integral, and $\cI\subset\cO_T^\gp$ is a
(quasi-coherent) fractional ideal of $\cO_T$, then $\cI$ is
coherent, provided the stalks $\cI_t$ are finitely generated
$\cO_{T,t}$-modules, for every $t\in T$. (Details left to the
reader.)
\end{remark}

\begin{remark}\label{rem_div-on-fans}
(i)\ \
Let $T$ be any integral fan. We define a sheaf $\cDiv_T$ on $T$,
by letting $\cDiv_T(U)$ be the set of all reflexive fraction
ideals of $\cO_{\!U}$, for every open subset $U\subset T$.

(ii)\ \
Now, suppose that $T$ is locally fine; in this case, we can
endow $\cDiv_T$ with a natural structure of $T$-monoid, as
follows. First, we define a presheaf of monoids on the site
$\cC_T$ of affine open subsets of $T$ by the rule :
$$
U\mapsto\cDiv_T(U):=(\Div(\cO_T(U)),\odot)
$$
(notation of \eqref{subsec_Div}) and for an inclusion $U'\subset U$
of affine open subset, the corresponding morphism of monoids
$\cDiv_T(U)\to\cDiv_T(U')$ is deduced from the flat map
$\cO_T(U)\to\cO_T(U')$, by virtue of lemma \ref{lem_rflx-rflx}(iv).
Arguing as in \eqref{subsec_def-Spec-P}, we see that $\cDiv_T$ is
a sheaf on $\cC_T$, and then \cite[Ch.0, \S3.2.2]{EGAI} implies
that $\cDiv_T$ extends uniquely to a sheaf of monoids on $T$.
It is then clear that the sheaf of sets underlying this $T$-monoid
is (naturally isomorphic to) the sheaf defined in (i).

(iii)\ \
In the situation of (ii), we have likewise a $T$-submonoid
$\cDiv_T^+\subset\cDiv_T$ (remark \ref{rem_pretty-obvious}),
and we may also define a $T$-monoid $\bar\cDiv_T$ (see
\eqref{subsec_J_P-plus}). Moreover, we have the global version
of \eqref{eq_basic-Div-seq} : namely, the sequence of $T$-monoids
$$
1\to\cO_T^\gp\xrightarrow{\ j_T\ }\cDiv_T\to
\bar\cDiv_T\to 1
$$
is exact (recall that $\cO_T^\times=1_T$, the initial $T$-monoid),
and $j_T$ restricts to a map of $T$-monoids
$$
\cO_T\to\cDiv^+_T.
$$
Indeed, the assertion can be checked on the stalks
over each $t\in T$, where it reduces to the exact sequence
\eqref{eq_basic-Div-seq} for $P:=\cO_{T,t}$. Lastly, we
remark that if $T$ is locally fine and saturated, then
$\cDiv_T$ and $\bar\cDiv_T$ are abelian $T$-groups
(proposition \ref{prop_completely-sat}(i,ii)).
\end{remark}

\sset\subsubsection{}
If $f:T'\to T$ is a morphism of fans, and $\cM$ is any
$\cO_T$-module, then we define as usual the $\cO_{T'}$-module :
$$
f^*\!\cM:=f^{-1}\!\cM\otimes_{f^{-1}\cO_T}\cO_{T'}
$$
where $f^{-1}\!\cM$ denotes the usual sheaf-theoretic inverse
image of $\cM$ (so $f^{-1}\cO_T$ means here what was denoted
$f^*\cO_T$ in definition \ref{def_mon-spaces}(ii)). The rule
$\cM\mapsto f^*\!\cM$ yields a left adjoint to the functor
$$
\cO_{T'}\Mod\to\cO_T\Mod
\qquad
\cN\mapsto f_*\cN
$$
(verification left to the reader). Notice that if $\cM$
is quasi-coherent, then $f^*\!\cM$ is a quasi-coherent
$\cO_{T'}$-module. Indeed, the assertion is local on $T'$,
hence we are reduced to the case where $T'=(\Spec\,P')$
and $T=(\Spec\,P)$ for some monoids $P$ and $P'$. In this
case, the functor $M\mapsto f^*(M^\sim):P\Mod\to\cO_{T'}\Mod$
is left adjoint to the functor $\cM\mapsto\Gamma(T',\cM)$
on $\cO_{T'}$-modules. The latter functor also admits
the left adjoint given by the rule :
$M\mapsto(M\otimes_PP')^\sim$, whence a natural isomorphism
of $\cO_{T'}$-modules :
$$
f^*(M^\sim)\isom(M\otimes_PP')^\sim.
$$

\sset\subsubsection{}\label{subsec_height-in-T}
Let $T:=(T,\cO_T)$ be a fan, $t\in T$ any point. The {\em height\/}
of $t$ is :
$$
\hgt_T(t):=\dim\cO_{T,t}\in\N\cup\{+\infty\}
$$
(see definition \ref{def_dim-mon}) and the {\em dimension\/} of $T$
is $\dim T:=\sup(\hgt_T(t)~|~t\in T)$. (If $T=\emptyset$ is the
empty fan, we let $\dim T:=-\infty$.)

Suppose that $T$ is locally finite; then it follows from
\eqref{eq_ident-stalks} and lemma \ref{lem_face}(iii),(iv) that
the height of any point of $T$ is an integer. Moreover, let
$U(t)\subset T$ denote the subset of all points $x\in T$ which
specialize to $t$ ({\em i.e.} such that the topological closure
of $\{x\}$ in $T$ contains $t$); clearly $U(t)$ is the intersection
of all the open neighborhoods of $t$ in $T$, and we have a natural
homeomorphism :
\set\begin{equation}\label{eq_homeom-fans}
\Spec\,\cO_{T,t}\isom U(t).
\end{equation}
Especially, if $T$ is locally finite, $U(t)$ is a finite set, and
moreover $U(t)$ is an open subset : indeed, if $U\subset T$ is any
finite affine open neighborhood of $t$, we have $U(t)\subset U$,
hence $U(t)$ can be realized as the intersection of the finitely
many open neighborhoods of $t$ in $U$. In this case,
\eqref{eq_homeom-fans} induces an isomorphism of fans :
\set\begin{equation}\label{eq_extisom-fans}
(\Spec\,\cO_{T,t})^\sharp\isom(U(t),\cO_{T|U(t)}).
\end{equation}
Therefore, for every $h\in\N$, let $T_h\subset T$ be the subset of
all points of $T$ of height $\leq h$; clearly $U(t)\subset T_h$
whenever $t\in T_h$, hence the foregoing shows that -- if $T$ is
locally finite -- $T_h$ is an open subset of $T$ for every $h\in\N$,
and $T=\bigcup_{h\in\N}T_h$.

Notice also that the simplicial locus of a fan $T$ is closed under
generizations. Therefore, $T_\mathrm{sim}$ is an open subset of
$T$, whenever $T$ is locally finite.

\sset\subsubsection{}\label{prod_of-fine-fans} In the situation
of remark \ref{rem_product-of-aff-fans}(v), suppose additionally
that $P$ and $Q$ are finitely generated. The natural projection
$P\times Q\to P$ induces a morphism
$j:T_P\to T_{P\times Q}$ of affine fans, and it is easily
seen that $j(t)=\pi^{-1}(t,\emptyset)$ for every $t\in T_P$.
It follows that the restriction of $j$ is a homeomorphism
$U(t)\isom U(j(t))$ for every $t\in T_P$, and moreover
$\log j:j^*\cO_{T_{P\times Q}}\to\cO_{T_P}$ is an isomorphism.
We conclude that $j$ is an open immersion.

\sset\subsubsection{}\label{subsec_from-face-to-prime}
Let $P$ be a fine, sharp and saturated monoid, and set $Q:=P^\vee$
(notation of \eqref{subsec_dual-of-mon}). By proposition
\ref{prop_reflex-dual}(iv), we have a natural identification
$P\isom Q^\vee$. By proposition \ref{prop_Gordie}(ii) and corollary
\ref{cor_spanning-cone}(ii), the rule
\set\begin{equation}\label{eq_face-to-face}
F\mapsto F^*:=F_\R^*\cap P
\end{equation}
establishes a natural bijection from the faces of $Q$ to those of
$P$. For every face $F$ of $Q$, set $\fp_F:=P\!\setminus\!F^*$; there
follows a natural bijection $F\mapsto\fp_F$ between the set of all
faces of $Q$ and $\Spec\,P$, such that
$$
F\subset F'\Leftrightarrow\fp_F\subset\fp_{F'}.
$$
Moreover, set $T_P:=(\Spec\,P)^\sharp$; we have natural
identifications :
$$
F^\vee\isom\cO_{T_P,\fp_F} \qquad \text{for every face $F$ of $Q$}
$$
under which, the specialization maps
$\cO_{T_P,\fp_{F'}}\to\cO_{T_P,\fp_F}$ correspond to the restriction
maps $(F')^\vee\to F^\vee$ : $\phi\mapsto\phi_{|F}$.

\begin{definition}\label{def_subdivisions}
Let $T:=(T,\cO_T)$ be a fan.
\begin{enumerate}
\item
An {\em integral} (resp. a {\em rational}) {\em partial subdivision}
of $T$ is a morphism $f:(T',\cO_{T'})\to T$ of fans such that, for
every $t\in T'$, the group homomorphism
$$
(\log f)^\gp_t:\cO^\gp_{T,f(t)}\to\cO^\gp_{T',t} \qquad \text{(resp.
the $\Q$-linear map \ \ $(\log f)^\gp_t\otimes_\Z\one_\Q$)}
$$
is surjective.
\item
If $f:T'\to T$ is an integral (resp. rational) partial subdivision,
and the induced map
$$
T'(\N)\to T(\N) \qquad \text{(resp. $T'(\Q_+)\to T(\Q_+)$)} \quad :
\quad \phi\mapsto f\circ\phi
$$
is bijective, we say that $f$ is an integral (resp. a rational)
{\em subdivision} of $T$.
\item
A morphism of fans $f:T'\to T$ is {\em finite} (resp. {\em proper}),
if the fibre $f^{-1}(t)$ is a finite (resp. and non-empty) set, for
every $t\in T$.
\item
A subdivision $T'\to T$ of\/ $T$ is {\em simplicial}, if
$T'_\mathrm{sim}=T'$.
\end{enumerate}
\end{definition}

\begin{remark}(i)\ \ Let $T$ be any integral fan. Then the counit
of adjunction
$$
T^\sat\to T
$$
is an integral subdivision. This morphism is also a homeomorphism
on the underlying topological spaces, in light of
\ref{lem_exc-satura}(iv).

(ii)\ \ Let $f:T'\to T$ be any integral subdivision. Then $f$ restricts
to a bijection $T'_0\isom T_0$ on the sets of points of height zero.
Indeed, notice that if $t'\in T'$ is a point of height zero, then
$f(t')\in T_0$ since the map $(\log f)^\gp_{t'}$ must be local;
moreover $\mathrm{loc}.\Hom_\Mnd(\cO_{T,t'},\N)$ consists of precisely
one element, namely the unique map $\sigma_t:\N\to\{1\}$, and if
$t'_1,t'_2\in T'_0$ have the same image in $T_0$, the sections
$\sigma_{t'_1}$ and $\sigma_{t'_2}$ have the same image in $T(\N)$,
hence they must coincide, so that $t'_1=t'_2$, as claimed.

(iii)\ \ Let $f:T'\to T$ be an integral subdivision of locally fine
and saturated fans. In general, the image of a point $t'\in T'$ of
height one may have height strictly greater than one. On the other
hand, for any $t\in T$ of height one, and any $t'\in f^{-1}(t)$, the
map $\Z\isom\cO_{T,t}^\gp\to\cO_{T',t'}^\gp$ must be surjective
(theorem \ref{th_structure-of-satu}(ii)), therefore $\cO_{T',t'}^\gp$
is a cyclic group; however $\cO_{T',t'}$ is also sharp and saturated,
so it must be either the trivial monoid $\{1\}$ or $\N$.
The first case is excluded by (ii), so $\hgt(t')=1$, and moreover
$(\log f)_{t'}$ is an isomorphism (and there exists a unique such
isomorphism). Since the induced map $T'(\N)\to T(\N)$ is bijective,
it follows easily that $f^{-1}(t)$ consists of exactly one point, and
therefore $f$ restricts to an isomorphism $f^{-1}(T_1)\isom T_1$.
\end{remark}

\begin{proposition}\label{prop_val-crit-fans}
Let $f:T'\to T$ be a morphism of fans, with $T'$ locally finite,
and consider the following conditions :
\begin{enumerate}
\alphaenu
\item
The induced map $T'(\N)\to T(\N)$ is injective.
\item
For every integral saturated monoid $P$, the induced map
$T'(P)\to T(P)$ is injective.
\item
$f$ is a partial rational subdivision.
\end{enumerate}
Then we have : {\em (a) $\Leftrightarrow$ (b) $\Rightarrow$ (c)}.
\end{proposition}
\begin{proof} Obviously (b) $\Rightarrow$ (a). Conversely, assume
that (a) holds, let $P$ be a saturated monoid, and suppose we have
two sections in $T'(P)$ whose images in $T(P)$ agree. In light of
example \ref{ex_fans}(iii), this means that we may find two points
$t'_1,t'_2\in T'$, such that $f(t'_1)=f(t_2')=t$, and two local
morphisms of monoids $\sigma_i:\cO_{T',t'_1}\to P/P^\times$ whose
compositions with $\log f_{t'_i}$ ($i=1,2$) yield the same morphism
$\cO_{T,t}\to P/P^\times$, and we have to show that these maps are
equal.
In view of lemma \ref{lem_exc-satura}(ii), we may then replace $P$
by $P/P^\times$, and assume that $P$ is sharp.
Since the stalks of $\cO_{\!T'}$ are finitely generated, the morphisms
$\sigma_i$ factor through a finitely generated submonoid $M\subset P$.
We may then replace $P$ by its submonoid $M^\sat$, which allows to
assume additionally that $P$ is finitely generated (corollary
\ref{cor_fragment-Gordon}(ii)). In this case, we may find an injective
map $j:P\to\N^{\oplus r}$ (corollary \ref{cor_consequent}(iv); notice
that $j$ is trivially a local morphism), hence we may replace
$\sigma_i$ by $j\circ\sigma_i$ (for $i=1,2$), after which we may
assume that $P=\N^{\oplus r}$ for some $r\in\N$. Let $\delta:P\to\N$
be the local morphism given by the rule :
$(x_1,\dots,x_r)\mapsto x_1+\dots+x_r$ for every $x_1,\dots,x_r\in\N$;
the compositions $\delta\circ\sigma_i$ (for $i=1,2$) are two elements
of $T'(\N)$ whose images agree in $T(\N)$, hence they must coincide
by assumption. This implies already that $t'_1=t'_2$. Next, let
$\pi_k:P\to\N$ (for $k=1,\dots,r$) be the natural projections, and
fix $k\leq r$; the morphisms $\pi_k\circ\sigma_i$ for $i=1,2$ are
not necessarily local, but they determine elements of $T'(\N)$ whose
images agree again in $T(\N)$, hence they must coincide. Since $k$ is
arbitrary, we deduce that $\sigma_1=\sigma_2$, as stated.

Next, we suppose that (b) holds, and we wish to show assertion (c);
the latter is local on $F'$, hence we may assume that both $F$ and
$F'$ are affine, say $F=(\Spec\,Q)^\sharp$ and $F'=(\Spec\,Q')^\sharp$,
with $Q'$ finitely generated and sharp, and then we are reduced to
checking that the map
$Q^\gp\otimes_\Z\Q\to Q^{\prime\gp}\otimes_\Z\Q$ induced by $f$
is surjective, or equivalently, that the dual map :
$$
\Hom_\Mnd(Q',\Q)\to\Hom_\Mnd(Q,\Q)
$$
is injective. To this aim, we may further assume that $Q'$ is
integral, in which case, by remark \ref{rem_reflex-dual}(i) we have
$\Hom_\Mnd(Q',\Q)=\Hom_\Mnd(Q',\Q_+)^\gp$; the contention is an easy
consequence.
\end{proof}

\sset\subsubsection{}
In light of proposition \ref{prop_val-crit-fans}, we may
ask whether the surjectivity of the map on $P$-points induced
by a morphism $f$ of fans can be similarly characterized.
This turns out to be the case, but more assumptions must be
made on the morphism $f$, and also some additional
restrictions must be imposed on the type of monoid $P$.
Namely, we shall consider monoids of the form $\Gamma_{\!+}$,
where $(\Gamma,\leq)$ is any totally ordered abelian group,
and $\Gamma_{\!+}\subset\Gamma$ is the subgroup of all elements
$\leq 1$ (where $1\in\Gamma$ denotes the neutral element).
With this notation, we have the following :

\begin{proposition}\label{prop_order-only}
Let $f:T'\to T$ be a finite partial integral subdivision, with $T$
locally finite. The following conditions are equivalent :
\begin{enumerate}
\alphaenu
\item
The induced map $T'(\N)\to T(\N)$ is surjective.
\item
For every totally ordered abelian group $(\Gamma,\geq)$, the
induced map $T'(\Gamma_{\!+})\to T(\Gamma_{\!+})$ is surjective.
\end{enumerate}
\end{proposition}
\begin{proof} Obviously, we need only to show that (a)
$\Rightarrow$ (b). Thus suppose, by way of contradiction, that
(a) holds, but nevertheless there exist a totally ordered
abelian group $(\Gamma,\leq)$, and an element of $T(\Gamma_{\!+})$
which is not in the image of $T'(\Gamma_{\!+})$. Such element
corresponds to a local morphism of monoids
$\phi:\cO_{T,t}\to\Gamma_{\!+}$, for some $t\in T$, and the
assumption means that $\phi$ does not factor through the monoid
$\cO_{T',s}$, for any $s\in f^{-1}(t)$. Set
$$
P:=\cO_{T,t}^\intg \qquad
Q_s:=P^\gp\times_{\cO^\gp_{\!T',s}}\cO^\intg_{T',s}
\quad\text{for every $s\in f^{-1}(t)$}
$$
Notice that, since the map $P^\gp\to\cO^\gp_{\!T',s}$ is surjective,
we have
$$
Q_s/Q_s^\times\simeq\cO^\intg_{T',s}/(\cO^\intg_{\!T',s})^\times
$$
and there is a natural injective morphism of monoids $g_s:P\to Q_s$,
determined by the pair $(i,(\log f)^\intg_s)$, where $i:P^\intg\to
P^\gp$ is the natural morphism; moreover, $P^\gp=Q_s^\gp$ for every
$s\in f^{-1}(t)$. Clearly $\phi$ factors through a morphism
$\bar\phi:P\to\Gamma_{\!+}$; since the unit of adjunction
$\cO_{T,t}\to P$ is surjective, it follows that $P$ is sharp and
$\bar\phi$ is local. Moreover, the group homomorphism
$\bar\phi{}^\gp:P^\gp\to\Gamma$ factors uniquely through each
$Q_s^\gp$. Our assumption then states that we may find, for each
$s\in f^{-1}(t)$, an element $x_s\in Q_s$ whose image in $\Gamma$
lies in the complement of $\Gamma_{\!+}$, {\em i.e.} the image of
$x_s^{-1}$ lies in the maximal ideal $\fm\subset\Gamma_{\!+}$. Let
$P'\subset P^\gp$ be the submonoid generated by $P$ and by
$(x_s^{-1}~|~s\in f^{-1}(t))$. By construction, $P'$ is finitely
generated, and the morphism $\bar\phi$ extends uniquely to a
morphism $P'\to\Gamma_{\!+}$, which maps each $x_s^{-1}$ into $\fm$.
It follows that all the $x_s^{-1}$ lie in the maximal ideal of $P'$.
Let us now pick any local morphism $\psi':P'\to\N$ (corollary
\ref{cor_consequent}(iii)); by restriction, $\psi'$ induces a local
morphism $\psi:P\to\N$, which -- according to (a) -- must factor
through a local morphism $\psi_s:Q_s\to\N$, for at least one $s\in
f^{-1}(t)$. However, on the one hand we have
$\psi^{\prime\gp}=\psi^\gp=\psi^\gp_s$; on the other hand
$\psi'(x_s^{-1})\neq 0$, hence
$\psi_s(x_s)=\psi^{\prime\gp}(x_s)\notin\N$, a contradiction.
\end{proof}

\begin{example}\label{ex_uniqueness-of-mu_n}
Let $T$ be a locally fine fan, $\phi:F\to T$ an integral subdivision,
with $F$ locally fine and saturated, and $k>0$ an integer. Suppose we
have a commutative diagram of fans :
\set\begin{equation}\label{eq_yad-of-fans}
{\diagram  F \ar[r]^-{\phi} \ar[d]_g & T \ar[d]^{\bek_T} \\
                   F \ar[r]^-{\phi} & T.
\enddiagram}
\end{equation}
where $\bek_T$ is the $k$-Frobenius endomorphism
(example \ref{ex_multiply-by-n-in-fan}(ii)). Then we claim
that necessarily $g=\bek_F$. Indeed, suppose that this fails;
then we may find a point $t\in F$ such that the composition
of $g$ and the open immersion $j_t:U(t)\to F$ is not equal to 
$j_t\circ\bek_{U(t)}$. Set $P:=\cO_{F,t}$; then
$g\circ j_t\neq j_t\circ\bek_{U(t)}$ in $F(P)$. However, an
easy  computation shows that
$\phi(g\circ j_t)=\phi(j_t\circ\bek_{U(t)})$, which
contradicts proposition \ref{prop_val-crit-fans}.
\end{example}

\sset\subsubsection{}\label{subsec_proj-fan}
Let $P:=\coprod_{n\in\N}P_n$ be a $\N$-graded monoid (see
definition \ref{def_grad-monoids}); then $P_0$ is a submonoid
of $P$, every $P_n$ is a $P_0$-module, and $P_+:=\coprod_{n>0}P_n$
is an ideal of $P$. For every $a\in P$, the localization $P_a$ is
$\Z$-graded in an obvious way, and we denote by $P_{(a)}\subset P_a$
the submonoid of elements of degree $0$. Notice that there is a
natural identification of $P_0$-monoids
\set\begin{equation}\label{eq_invar-by-pows}
P_{(a^n)}\isom P_{(a)}
\qquad
\text{for every integer $n>0$.}
\end{equation}
Set as well :
$$
D_+(a):=(\Spec\,P_{(a)})^\sharp
$$
and notice that the natural map $P\to P_{(a)}$ induces a
morphism of fans $\pi_a:D_+(a)\to T_P:=(\Spec\,P_0)^\sharp$. If
$b\in P$ is any other element, in order to determine the fibre
product $D_+(a)\times_{T_P}D_+(b)$ we may assume -- in light
of \eqref{eq_invar-by-pows} -- that $a,b\in P_n$ for the same
integer $n$, in which case we have natural isomorphisms
$$
P_{(a)}\otimes_{P_0}P_{(b)}\isom P_{(ab)}\xleftarrow{\sim}
P_{(a)}[b^{-1}a]
$$
(see remark \ref{rem_tens-is-pushout}(i)) onto the localization of
$P_{(a)}$ obtained by inverting its element $a^{-1}b$; this is
of course the same as $P_{(b)}[a^{-1}b]$. In other words
$D_+(a)\times_{T_P}D_+(b)$ is naturally isomorphic to $D_+(ab)$,
identified to an open subfan in both $D_+(a)$ and $D_+(b)$.
We may then glue the fans $D_+(a)$ for $a$ ranging over all
the elements of $P$, to obtain a new fan, denoted :
$$
\Proj\,P
$$
called the {\em projective fan\/} associated with $P$.
By inspecting the construction, we see that the morphisms
$\pi_a$ assemble to a well defined morphism of fans
$\pi_P:\Proj\,P\to T_P$. Each element $a\in P$ yields an open
immersion $j_a:D_+(a)\to\Proj\,P$, and  if $b\in P$ is any
other element, $j_{ab}$ factors through an open immersion
$D_+(ab)\to D_+(a)$.

\sset\subsubsection{}\label{subsec_morph-proj-fan}
Let $\phi:P\to P'$ be a morphism of $\N$-graded monoids (so
$\phi P_n\subset P'_n$ for every $n\in\N$). Set :
$$
G(\phi):=\bigcup_{a\in P}D_+(\phi(a))\subset\Proj\,P'.
$$
Notice that, for every $a\in P$, $\phi$ induces a morphism
$\phi_{(a)}:P_{(a)}\to P'_{(a)}$, whence a morphism of
affine fans $(\Proj\,\phi)_a:D_+(\phi(a))\to D_+(a)\subset\Proj\,P$.
Moreover, if $b\in P$ is any other element, it is easily
seen that $(\Proj\,\phi)_a$ and $(\Proj\,\phi)_b$ agree
on $D_+(\phi(a))\cap D_+(\phi(b))$. Therefore, the morphisms
$(\Proj\,\phi)_a$ glue to a well defined morphism :
$$
\Proj\,\phi:G(\phi)\to\Proj\,P.
$$
Notice that $G(\phi)=\Proj\,P'$, whenever $\phi P$ generates
the ideal $P'_+$. Moreover, we have
\set\begin{equation}\label{eq_pull-back-basic}
(\Proj\,\phi)^{-1}D_+(a)=D_+(\phi(a))
\qquad\text{for every $a\in P$.}
\end{equation}
Indeed, say that $D_+(b)\subset G(\phi)$ for some $b\in P'$,
and $(\Proj\,\phi)D_+(b)\subset D_+(a)$. In order to show that
$D_+(b)\subset D(\phi(a))$, it suffices to check that
$D_+(b\phi(c))\subset D_+(\phi(a))$ for every $c\in P$.
However, the assumption means that the natural map
$$
P_{(c)}\to P'_{(\phi(c))}\to P'_{(b\phi(c))}\to
P'_{(b\phi(c))}/P^{\prime\times}_{(b\phi(c))}
$$
factors through the localization $P_{(c)}\to P_{(ac)}$. This
is equivalent to saying that $\phi(c^{-1}a)$ is invertible
in $P'_{(b\phi(c))}$, in which case the localization
$P'_{(\phi(c))}\to P'_{(b\phi(c))}$ factors through
the localization $P'_{(\phi(c))}\to P'_{(\phi(ac))}$.
The latter means that the open immersion
$D_+(b\phi(c))\subset D_+(\phi(c))$ factors through the
open immersion $D_+(\phi(ac))\subset D_+(\phi(c))$, as
claimed.

\sset\subsubsection{}\label{subsec_notate-U_n}
In the situation of \eqref{subsec_proj-fan}, set $Y:=\Proj\,P$
to ease notation. Let $M$ be a $\Z$-graded $P$-module; for every
$a\in P$, let $M_{(a)}\subset M_a:=M\otimes_PP_a$ be the
$P_{(a)}$-submodule of degree zero elements (for the natural
grading on $M_a$).
We deduce a quasi-coherent $\cO_{D_+(a)}$-module $M^\sim_{(a)}$
(see definition \ref{def_quasicoh-fans}). Moreover, if $b\in P$
is any other element, we have a natural identification
$$
\tilde\omega_{a,b}:
M^\sim_{(a)|D_+(a)\cap D_+(b)}\isom M^\sim_{(b)|D_+(a)\cap D_+(b)}.
$$
This can be verified as follows. First, in view of
\eqref{eq_invar-by-pows}, we may assume that $a,b\in P_n$,
for some $n\in\N$, in which case we consider the $P$-linear
morphism :
$$
M_{(a)}\to M_{(b)}\otimes_{P_{(b)}}P_{(b)}[a^{-1}b]
\quad :\quad
\frac{x}{a^m}\mapsto\frac{x}{b^m}\otimes\frac{b^m}{a^m}
\qquad\text{for every $x\in M_{nm}$}.
$$
It is easily seen that this map is actually $P_{(a)}$-linear,
hence it extends to a $\cO_T(D_+(a)\cap D_+(b))$-linear morphism :
$$
\omega_{a,b}: M_{(a)}\otimes_{P_{(a)}}P_{(a)}[b^{-1}a]\isom
 M_{(b)}\otimes_{P_{(b)}}P_{(b)}[a^{-1}b].
$$
Moreover, $\omega_{a,b}\circ\omega_{b,a}$ is the identity
map, hence $\omega_{a,b}$ induces the sought isomorphism
$\tilde\omega_{a,b}$. Furthermore, for any $a,b,c\in P$, set
$D_+(a,b,c):=D_+(a)\cap D_+(b)\cap D_+(c)$; we have the identity :
$$
\tilde\omega_{a,c|D_+(a,b,c)}=
\tilde\omega_{b,c|D_+(a,b,c)}\circ\tilde\omega_{a,b|D_+(a,b,c)}
$$
which shows that the locally defined sheaves $M^\sim_{(a)}$
glue to a well defined $\cO_Y$-module, which we shall
denote $M^\sim$. Especially, for every $n\in\Z$, let $P(n)$ be
the $\Z$-graded $P$-module such that $P(n)_k:=P_{n+k}$ for
every $k\in\Z$ (with the convention that $P_n:=\emptyset$
if $n<0$); we set :
$$
\cO_Y(n):=P(n)^\sim.
$$
Every element $a\in P_n$ induces a natural isomorphism :
$$
\cO_Y(n)_{|D_+(a)}\isom\cO_{D_+(a)}
\quad :\quad x\mapsto f^{-k}x \quad\text{for every local section $x$.}
$$
Hence on the open subset :
$$
U_n(P):=\bigcup_{a\in P_n}D_+(a)
$$
the sheaf $\cO_Y(n)$ restricts to an invertible $\cO_{U_n(P)}$-module
(see definition \ref{def_coh-idea-log}(iv)). Especially, if $P_1$
generates $P_+$, the $\cO_Y$-modules $\cO_Y(n)$ are invertible,
for every $n\in\Z$.

\sset\subsubsection{}\label{subsec_sometimes-isom}
In the situation of \eqref{subsec_morph-proj-fan}, let $M$ be
a $\Z$-graded $P$-module. Then $M':=M\otimes_PP'$ is a $\Z$-graded
$P'$-module, with the grading defined by the rule :
\set\begin{equation}\label{eq_def-gradings}
M'_n:=\bigcup_{j+k=n}\Img(M_j\otimes_{P_0}P'_k\to M').
\end{equation}
There follows a $P_{(a)}$-linear morphism :
\set\begin{equation}\label{eq_pulbak}
M_{(a)}\to M'_{(\phi(a))}
\quad : \quad
\frac{x}{a^k}\mapsto\frac{x\otimes 1}{\phi(a)^k}
\qquad\text{for every $a\in P$}
\end{equation}
and since both localization and tensor product commute
with arbitrary colimits, it is easily seen that \eqref{eq_pulbak}
extends an injective $P'_{(\phi(a))}$-linear map
$$
M_{(a)}\otimes_{P_{(a)}}P'_{(\phi(a))}\to M'_{(\phi(a))}
$$
whence a map of $\cO_{D_+(\phi(a))}$-modules
$(\Proj\,\phi)^*M^\sim_{|D_+(\phi(a))}\to(M')^\sim_{|_{D_+(\phi(a))}}$,
and the system of such maps, for $a$ ranging over the elements
of $P$, is compatible with all open immersions
$D_+(\phi(ab))\subset D_+(a)$, whence a well defined monomorphism
of $\cO_{G(\phi)}$-modules
\set\begin{equation}\label{eq_global-pulbak}
(\Proj\,\phi)^*M^\sim\to(M')^\sim_{|G(\phi)}.
\end{equation}
Moreover, if $a\in P_1$, then for every $m\in M_j$ and $x\in P'_k$, we
may write
$$
\frac{m\otimes x}{\phi(a)^{j+k}}=\frac{m}{a^j}\otimes\frac{x}{\phi(a)^j}
$$
so the above map is an isomorphism on $D_+(a)$. Thus,
\eqref{eq_global-pulbak} restricts to an isomorphism on
the open subset
$$
G_1(\phi):=\bigcup_{a\in P_1}D_+(\phi(a)).
$$
Especially \eqref{eq_global-pulbak} is an isomorphism whenever
$P_1$ generates $P_+$. Notice as well that
$G_1(\phi)\subset U_1(P')\cap G(\phi)$, and actually
$G_1(\phi)=U_1(P')$ if $\phi(P)$ generates $P'_+$.

\sset\subsubsection{}\label{subsec_proj-and-basechange}
Let $P$ be as in \eqref{subsec_proj-fan}, and $f:P_0\to Q$
a given morphism of monoids. Then $P':=P\otimes_{P_0}Q$ is
naturally $\N$-graded, so that the natural map $f_P:P\to P'$
is a morphism of graded monoids. Every element of $P'$
is of the form $a\otimes b=(a\otimes 1)\cdot(1\otimes b)$,
where $a\in P$ and $b\in Q$. Then lemma \ref{lem_localize}
yields a natural isomorphism of $Q$-monoids :
$$
P_{(a)}\otimes_{P_0}Q[b^{-1}]\isom P'_{(a\otimes b)}
$$
whence an isomorphism of affine fans :
$$
\beta_{a\otimes b}:D_+(a\otimes b)\isom D_+(a)\times_{T_P}D(b)
$$
such that
$(\pi_a\times_{T_P}j^*_b)\circ\beta_{a\otimes b}=\pi_{a\otimes b}$,
where $j^*_b:D(b)\to T_Q:=(\Spec\,Q)^\sharp$ is the natural open
immersion. Especially, it is easily seen that the isomorphisms
$\beta_{a\otimes 1}$ assemble to a well defined isomorphism of fans :
$$
(\Proj\,f_P,\pi_{P'}):\Proj\,P'\isom\Proj\,P\times_{T_P}T_Q
$$
such that
$(\pi_P\times_{T_P}\one_{T_Q})\circ(\Proj\,f_P,\pi_{P'})=\pi_{P'}$.
Lastly, if $g:Q\to R$ is another morphism of monoids,
$T_R:=(\Spec\,R)^\sharp$, and $P'':=P'\otimes_QR$, then we have
the identity :
\set\begin{equation}\label{eq_compat-proj-fan}
((\Proj\,f_P,\pi_{P'})\times_{T_Q}\one_{T_R})\circ
(\Proj\,g_{P'},\pi_{P''})=(\Proj\,(g\circ f)_P,\pi_{P''}).
\end{equation}
Moreover, for every $\Z$-graded module $M$, the map
$(\Proj\,f_P)^*M^\sim\to(M\otimes_{P_0}Q)^\sim_{|G(f_P)}$ of
\eqref{eq_global-pulbak} is an isomorphism, regardless of whether
or not $P_1$ generates $P_+$ (verification left to the reader).
Especially, we get a natural identification :
$$
(\Proj\,f_P)^*\cO_Y(n)\isom\cO_{Y'}(n)
\qquad
\text{for every $n\in\Z$}
$$
where $Y:=\Proj\,P$ and $Y':=\Proj\,P'$.

\sset\subsubsection{}
In the situation of \eqref{subsec_proj-and-basechange}, let
$\phi:R\to P$ be a morphism of $\N$-graded monoids. There
follow morphisms of fans :
$$
\Proj\,\phi:G(\phi)\to\Proj\,R
\qquad
\Proj(f_P\circ\phi):G(f_P\circ\phi)\to\Proj\,R
$$
and in view of \eqref{eq_pull-back-basic}, it is easily
seen that :
\set\begin{equation}\label{eq_courage-proj}
G(f_P\circ\phi)=(\Proj\,f_P)^{-1}(G(\phi)).
\end{equation}

\sset\subsubsection{}\label{subsec_general-proj}
Let now $(T,\cO_T)$ be any fan. A {\em $\N$-graded $\cO_T$-monoid}
is a $\N$-graded $T$-monoid $\cP$, with a morphism $\cO_T\to\cP$
of $T$-monoids. We say that such a $\cO_T$-monoid is
{\em quasi-coherent}, if it is such, when regarded as a
$\cO_T$-module. To a quasi-coherent $\N$-graded $\cO_T$-monoid
$\cP$, we attach a morphism of fans :
$$
\pi_\cP:\Proj\,\cP\to T
$$
constructed as follows. First, for every affine open subfan
$U\subset T$, the monoid $\cP(U)$ is $\N$-graded, so we have
the projective fan $\Proj\,\cP(U)$, and the morphism of monoids
$\cO_T(U)\to\cP(U)$ induces a morphism of fans
$\Proj\,\cP(U)\to U$. Next, say that $U_1,U_2\subset T$ are
two affine open subsets; for any affine open subset
$V\subset U_1\cap U_2$ we have restriction maps
$\rho_{V,i}:\cP(U_i)\to\cP(V)$ inducing isomorphisms of graded
$\cO_T(V)$-monoids :
$$
\cP(U_i)\otimes_{\cO_T(U_i)}\cO_T(V)\isom\cP(V).
$$
whence isomorphisms of $V$-fans :
$$
\Proj\,\cP(V)\isom\Proj\,\cP(U_i)\otimes_{\cO_T(U_i)}\cO_T(V)
\xrightarrow{ \Proj\,\rho_{V,i} }\Proj\,\cP(U_i)\times_{U_i}V
$$
which in turn yield natural identifications :
$$
\theta_V:\Proj\,\cP(U_1)\times_{U_1}V\isom\Proj\,\cP(U_2)\times_{U_2}V.
$$
If $W\subset V$ is a smaller affine open subset,
\eqref{eq_compat-proj-fan} implies that
$\theta_V\times_V\one_W=\theta_W$, and therefore the isomorphisms
$\theta_V$ glue to a single isomorphism of $U_1\cap U_2$-fans :
$$
\Proj\,\cP(U_1)\times_{U_1}(U_1\cap U_2)\isom
\Proj\,\cP(U_2)\times_{U_2}(U_1\cap U_2)
$$
which is furthermore compatible with base change to any
triple intersection $U_1\cap U_2\cap U_3$ of affine open
subsets (details left to the reader). In such situation,
we may glue the fans $\Proj\,\cP(U)$ -- with $U\subset T$
ranging over all the open affine subsets -- along the
above isomorphisms, to obtain the sought fan $\Proj\,\cP$;
the construction also comes with a well defined morphism
to $T$, as required. Then, for every such open affine $U$,
the induced morphism $\Proj\,\cP(U)\to\Proj\,\cP$ is an
open immersion; finally a direct inspection shows that,
for every smaller affine open subset $V\subset U$ we have :
$$
U_n(\cP(U))\cap\Proj\,\cP(V)=U_n(\cP(V))
\qquad\text{for every $n\in\N$}
$$
(where the intersection is taken in $\Proj\,\cP$). Hence the
union of all the open subsets $U_n(\cP(U))$ is an open subset
$U_n(\cP)\subset\Proj\cP$, intersecting each $\Proj\,\cP(U)$
in its subset $U_n(\cP(U))$.

\sset\subsubsection{}
To ease notation, set $Y:=\Proj\,\cP$, and let $\pi:Y\to T$
be the projection. Let $\cM$ be a $\Z$-graded
$\cP$-module, quasi-coherent as a $\cO_T$-module; for every
affine open subset $U\subset T$, the graded $\cP(U)$-module
$\cM(U)$ yields a quasi-coherent $\cO_{\pi^{-1}U}$-module
$\cM_U^\sim$, and every inclusion of affine open subset
$U'\subset U$ induces a natural isomorphism
$\cM^\sim_{U|U'}\isom\cM^\sim_{U'}$ of $\cO_{\pi^{-1}U'}$-modules.
Therefore the modules $\cM^\sim_U$ glue to a well defined
$\cO_Y$-module $\cM^\sim$.

For every $n\in\Z$, denote by $\cM(n)$ the $\Z$-graded
$\cP$-module such that $\cM(n)_k:=\cM_{n+k}$ for every $k\in\Z$
(especially, with the convention that $\cP_k:=0$ whenever $k<0$,
we obtain in this way the $\cP$-module $\cP(n)$). We set :
$$
\cO_Y(n):=\cP(n)^\sim
\quad\text{and}\quad
\cM^\sim(n):=\cM^\sim\otimes_{\cO_Y}\cO_Y(n).
$$
Clearly the restriction of $\cO_Y(n)$ to $U_n(\cP)$ is invertible,
for every $n\in\Z$.

Moreover, for every $n\in\Z$, the scalar multiplication
$\cP(n)\otimes_{\cO_T}\cM\to\cM(n)$ determines a well defined
morphism of $\cO_Y$-modules :
$$
\cM^\sim(n)\to\cM(n)^\sim
$$
and arguing as in \eqref{subsec_sometimes-isom} we see that
the restriction of this map is an isomorphism on $U_1(\cP)$.
Especially, we have natural morphisms of $\cO_Y$-modules :
$$
\cO_Y(n)\otimes_{\cO_Y}\cO_Y(m)\to\cO_Y(n+m)
\qquad
\text{for every $n,m\in\Z$}
$$
whose restrictions to $U_1(\cP)$ are isomorphisms.

\begin{example}\label{ex_symm-monoid}
Let $T$ be a fan, $\cL$ an invertible $\cO_T$-module, and
set $\cP(\cL):=\Sym^\bullet_{\cO_T}\cL$ (see example
\ref{ex_symm-pow}). Then the morphism
$$
\pi_{\cP(\cL)}:\P(\cL):=\Proj\,\cP(\cL)\to T
$$
is an isomorphism. Indeed, the assertion can be checked
locally on every affine open subset $U\subset T$, hence say
that $U=(\Spec\,P)^\sharp$ for some monoid $P$, and
$\cL\simeq\cO_{T|U}$, in which case the $P$-monoid $\cP(\cL)(U)$
is isomorphic to $P\times\N$ (with its natural morphism
$P\to P\times\N$ : $x\mapsto(x,0)$ for every $x\in P$),
and the sought isomorphism corresponds to the natural
identification :
\set\begin{equation}\label{eq_ident-trivial}
P=(P\times\N)_{(1,1)}
\end{equation}
where $(1,1)\in P\times\{1\}=(P\times\N)_1$. Likewise,
$\cO_{\P(\cL)}(n)$ is the $\cO_{\P(\cL)}$-module associated
with the graded $(P\times\N)$-module
$P\times\N(n)=\cL^{\otimes n}(U)\otimes_P(P\times\N)$, so
\eqref{eq_ident-trivial} induces a natural isomorphism
$$
\pi_{\cP(\cL)}^*\cL^{\otimes n}\isom\cO_{\P(\cL)}(n)
\qquad
\text{for every $n\in\N$.}
$$
\end{example}

\sset\subsubsection{}
In the situation of \eqref{subsec_general-proj}, let
$\phi:\cP\to\cP'$ be a morphism of quasi-coherent $\N$-graded
$\cO_T$-monoids (defined in the obvious way). By the foregoing,
for every affine open subset $U\subset T$, we have an induced
morphism $\Proj\,\phi(U):G(\phi(U))\to\Proj\cP(U)$ of $U$-fans,
where $G(\phi(U))\subset\Proj\,\cP(U)$ is an open subset of
$\Proj\,\cP'$. Let $V\subset U$ be a smaller affine open subset;
in light of \eqref{eq_courage-proj}, we have
$$
G(\phi(V))=G(\phi(U))\cap\Proj\,\cP(V).
$$
It follows that the union of all the open subsets $G(\phi(U))$
is an open subset $G(\phi)$ such that
$$
G(\phi)\cap\Proj\,\cP(U)=G(\phi(U))
\qquad
\text{for every affine open subset $U\subset T$}
$$
and the morphisms $\Proj\,\phi(U)$ assemble to a well defined
morphism
$$
\Proj\,\phi:G(\phi)\to\Proj\,\cP.
$$
Moreover, if $\cM$ is a $\Z$-graded quasi-coherent $\cP$-module,
the morphisms \eqref{eq_global-pulbak} assemble to a well defined
morphism of $\cO_{G(\phi)}$-modules :
\set\begin{equation}\label{eq_toberestricted}
(\Proj\,\phi)^*\cM^\sim\to(\cM')^\sim_{|G(\phi)}
\end{equation}
where the grading of $\cM':=\cM\otimes_\cP\cP'$ is defined as in
\eqref{eq_def-gradings}.
Likewise, the union of all open subsets $G_1(\phi(U))$ is an
open subset $G_1(\phi)\subset U_1(\cP)\cap G(\phi)$, such that
the restriction of \eqref{eq_toberestricted} to $G_1(\phi)$ is
an isomorphism. Especially, set $Y:=\Proj\,\cP$ and $Y':=\Proj\,\cP'$;
we have a natural morphism :
$$
(\Proj\,\phi)^*\cO_Y(n)\isom\cO_{Y'}(n)_{|G(\phi)}
$$
which is an isomorphism, if $\cP_1$ generates
$\cP_+:=\coprod_{n>0}\cP_n$ locally on $T$.

\sset\subsubsection{}\label{subsec_pull-back-proj-glob}
On the other hand, let $f:T'\to T$ be a morphism of fans.
The discussion in \eqref{subsec_proj-and-basechange} implies
that $f$ induces a natural isomorphism of $T'$-fans :
\set\begin{equation}\label{eq_burningout}
\Proj\,f^*\cP\isom\Proj\,\cP\times_TT'.
\end{equation}
Moreover, set $Y:=\Proj\,\cP$, $Y':=\Proj\,f^*\cP$, and let
$\pi_Y:Y'\to Y$ be the projection deduced from
\eqref{eq_burningout}; then there follows a natural identification :
$$
\cO_{Y'}(n)\isom\pi_Y^*\cO_Y(n)
\qquad\text{for every $n\in\Z$.}
$$

\sset\subsubsection{}\label{subsec_same-as-it-everwas}
Keep the notation of \eqref{subsec_general-proj}, and to ease
notation, set $Y:=\Proj\,\cP$.
Let $\cC$ be the category whose objects are all the pairs
$(\psi:X\to T,\cL)$, where $\psi$ is a morphism of fans, and
$\cL$ is an invertible $\cO_{\!X}$-module; the morphisms
$(\psi:X\to T,\cL)\to(\psi':X'\to T,\cL')$ are the pairs
$(\beta,h)$, where $\beta:X\to X'$ is a morphism of $T$-fans,
and $h:\beta^*\cL'\isom\cL$ is an isomorphism of
$\cO_{\!X'}$-modules (with composition of morphisms defined
in the obvious way). Consider the functor
$$
F_\cP:\cC^o\to\Set
$$
which assigns to any object $(\psi,\cL)$ of $\cC$, the set
consisting of all morphisms of graded $\cO_{\!X}$-monoids
$$
g:\psi^*\cP\to\Sym^\bullet_{\cO_{\!X}}\cL
$$
which are epimorphisms on the underlying $\cO_{\!X}$-modules
(notation of example \ref{ex_symm-pow}). On a morphism $(\beta,h)$
as in the foregoing, and an element $g'\in F_\cP(\psi',\cL')$,
the functor acts by the rule :
$$
F_\cP(\beta,h):=(\Sym^\bullet_{\cO_{\!X}}h)\circ\beta^*g'.
$$

\begin{lemma} In the situation of \eqref{subsec_same-as-it-everwas},
the following holds :
\begin{enumerate}
\item
The object
$(\pi_{\cP|U_1(\cP)}:U_1(\cP)\to T,\cO_Y(1)_{|U_1(\cP)})$
represents the functor $F_\cP$.
\item
If $\cP$ is an integral $T$-monoid, the $\cO_T$-monoid $\cP^\sat$
admits a unique grading such that the unit of adjunction
$\cP\to\cP^\sat$ is a $\N$-graded morphism, and there is a natural
isomorphism of\/ $\Proj\,\cP$-fans :
$$
\Proj\,(\cP^\sat)\isom(\Proj\,\cP)^\sat.
$$
\end{enumerate}
\end{lemma}
\begin{proof}(i): The proof is {\em mutatis mutandis}, the same
as that of lemma \ref{lem_repres-proj} (with some minor
simplifications). We leave it as an exercise for the reader.

(ii): The first assertion shall be left to the reader.
The second assertion is local on $\Proj\,\cP$, hence we may
assume that $T=(\Spec\,P_0)$, and $\cP=P^\sim$ for some $\N$-graded
integral $P_0$-monoid $P$. Let $a\in P^\sat$ be any element;
by definition we have $a^n\in P$ for some $n>0$, and we know
that the open subsets $D_+(a)$ et $D_+(a^n)$ coincide in
$\Proj(\cP^\sat)$; hence we come down to showing that
$(P_{(a)})^\sat=(P^\sat)_{(a)}$ for every $a\in P$, which
can be left to reader.
\end{proof}

\begin{definition} Let $(T,\cO_T)$ be a fan (resp. an integral fan),
$\cI\subset\cO_T$ an ideal (resp. a fractional ideal) of $\cO_T$.
\begin{enumerate}
\item
Let $f:X\to T$ be a morphism of fans (resp. of integral fans); then
$f^{-1}\!\cI$ is an ideal (resp. a fractional ideal) of $f^{-1}\cO_T$,
and we let :
$$
\cI\!\cO_{\!X}:=\log f(f^{-1}\!\cI)\cdot\cO_{\!X}
$$
which is the smallest ideal (res. fractional ideal) $\cO_{\!X}$
containing the image of $f^{-1}\cI$.
\item
A {\em blow up\/} of the ideal $\cI$ is a morphism of fans
(resp. of integral fans) $\phi:T'\to T$ which enjoys the
following universal property. The ideal (resp. fractional ideal)
$\cI\cO_{T'}$ is invertible, and every morphism of fans (resp.
of integral fans) $X\to T$ such that $\cI\cO_{\!X}$ is invertible,
factors uniquely through $\phi$.
\end{enumerate}
\end{definition}

\sset\subsubsection{}
Let $T$ be a fan (resp. an integral fan), $\cI\subset\cO_T$ a
quasi-coherent ideal (resp. fractional ideal), and consider the
$\N$-graded $\cO_T$-monoid :
$$
\cB(\!\cJ):=\coprod_{n\in\N}\cI^n
$$
where $\cI^n\subset\cO_T$ is the ideal (resp. fractional ideal)
associated with the presheaf $U\mapsto\cI(U)^n$ for every open
subset $U\subset T$ (notation of \eqref{subsec_toric}, with the
convention that $\cI^0:=\cO_T$) and the multiplication law of
$\cB(\cI)$ is defined in the obvious way.

\begin{proposition}
The natural projection
$$
\Proj\,\cB(\cI)\to T
$$
is a blow up of the ideal $\cI$.
\end{proposition}
\begin{proof} We shall consider the case where $T$ is not
necessarily integral, and $\cI\subset\cO_T$; the case of a
fractional ideal of an integral fan is proven in the same way.
Set $Y:=\Proj\,\cB(\cI)$; to begin with, let us show that
$\cI\!\cO_Y$ is invertible. The assertion is local on $T$,
hence we may assume that $T=(\Spec\,P)^\sharp$, and $\cI=I^\sim$
for some ideal $I\subset P$, so $Y=\Proj\,B(I)$, where
$B(I)=\coprod_{n\in\N}I^n$. Let $a\in B(I)_1=I$ be any
element; then the restriction of $\cI\!\cO_Y$ to $D_+(a)$ is
generated by $1=a/a\in B(I)_{(a)}$, so clearly
$\cI_{|D_+(a)}\simeq\cO_Y$; since $U_1(B(I))=Y$ (notation of
\eqref{subsec_notate-U_n}), the contention follows.

Next, let $\phi:X\to T$ be a morphism of fans, such that
$\cI\!\cO_{\!X}$ is an invertible ideal. It follows easily
that $\cI^n\cO_{\!X}$ is invertible for every $n\in\N$,
so the natural map of $\N$-graded $\cO_{\!X}$-monoids
$$
\phi^*\Sym^\bullet_{\cO_T}(\cI)\isom
\Sym^\bullet_{\cO_{\!X}}(\cI\!\cO_{\!X})\to\cB(\cI\!\cO_{\!X})
$$
is an isomorphism. On the other hand, the projection
$\Proj\,\cB(\cI\!\cO_{\!X})\to X$ is an isomorphism
(example \eqref{ex_symm-monoid}), whence -- in view of
\eqref{eq_burningout} -- a natural morphism of $T$-fans :
\set\begin{equation}\label{eq_only-you}
X\to\Proj\,\cB(\cI).
\end{equation}
To conclude, it remains to show that \eqref{eq_only-you}
is the only morphism of $T$-fans from $X$ to $\Proj\,\cB(\cI)$.
The latter assertion can be checked again locally on $T$,
so we are reduced as above to the case where $T$ is the spectrum
of $P$, and $\cI$ is associated with $I$. We may also assume
that $X=(\Spec\,Q)^\sharp$, and $\phi$ is given by a morphism
of monoids $f:P\to Q$. Then the hypothesis means that the
ideal $f(I)Q$ is isomorphic to $Q$ (see remark
\ref{rem_local-trivial-qcoh}(ii)), hence it is generated
by an element of the form $f(a)$, for some $a\in I$, and
the endomorphism $x\mapsto f(a)x$ of $f(I)Q$, is an isomorphism.
In such situation, it is clear that $f$ factors uniquely
through a morphism of monoids $P\to B(I)_{(a)}$; namely,
one defines $g:B(I)_{(a)}\to Q$ by the rule :
$a^{-k}x\mapsto f(a)^{-k}f(x)$ (for every $x\in I^k$), which
is well defined, by the foregoing observations. The morphism
$(\Spec\,g)^\sharp:X\to D_+(a)$ must then agree with
\eqref{eq_only-you}.
\end{proof}

\begin{example}\label{ex_blow-ups}
(i)\ \ 
Let $P$ be a monoid, $I\subset P$ any finitely generated
ideal, $\{a_1,\dots,a_n\}$ a finite system of generators of $I$; set
$T:=(\Spec\,P)^\sharp$, and let $\phi:T'\to T$ be the blow up of the
ideal $I^\sim\subset\cO_T$. Then $T'$ admits an open covering
consisting of the affine fans $D_+(f_i)$. The latter are the spectra
of the monoids $Q_i$ consisting of all fractions of the form $a\cdot
f_i^{-t}$, for every $a\in I^n$; we have $a\cdot f_i^{-t}=b\cdot
f_i^{-s}$ in $Q_i$ if and only if there exists $k\in\N$ such that
$a\cdot f_i^{s+k}=b\cdot f_i^{t+k}$, if and only if the two
fractions are equal in $P_{f_i}$, in other word, $Q_i$ is the
submonoid of $P_{f_i}$ generated by $P$ and $\{f_j\cdot
f_i^{-1}~|~j\leq n\}$, for every $i=1,\dots,n$.

(ii)\ \ 
Consider the special case where $P$ is fine, and the ideal
$I\subset P$ is generated by two elements $f,g\in P$. Let $t\in T'$
be any point; up to swapping $f$ and $g$, we may assume that $t$
corresponds to a prime ideal $\fp\subset P[f/g]$, hence $\phi(t)$
corresponds to $\fq:=j^{-1}\fp\subset P$, where $j:P\to P[f/g]$ is
the natural map. We have the following two possibilities :
\begin{itemize}
\item
Either $f/g\in\fp$, in which case let $y\in
F':=P[f/g]\!\setminus\!\fp$ be any element; writing
$x=y\cdot(f/g)^n$ for some $n\geq 0$ and $y\in P$, we deduce that
$n=0$, so $x=y\in F:=P\setminus\fq$, therefore $F'=j(F)$. Notice as
well that in this case $f/g$ is not invertible in $P[f/g]$, hence
$P[f/g]^\times=P^\times$, whence $\dim P[f/g]=\dim P$, and , by
corollary \ref{cor_consequent}(i).
\item
Or else $f/g\notin\fp$, in which case the same argument yields
$F'=j(F)[f/g]$. In this case, $f/g$ is invertible in $P[f/g]$ if and
only if it is invertible in the face $F'$, hence $\hgt\,\fp=\dim
P-\rk_\Z F'\geq\dim P-\dim F-1$, by corollary
\ref{cor_consequent}(i),(ii).
\end{itemize}
In either event, corollary \ref{cor_consequent}(i),(ii) implies the
inequality :
$$
1\geq\hgt(\phi(t))-\hgt(t)\geq 0 \qquad \text{for every $t\in T'$}.
$$
\end{example}

\subsection{Special subdivisions}\label{sec_special-sub}
In this section we explain how to construct -- either by
geometrical or combinatorial means -- useful subdivisions of
given fans.

\sset\subsubsection{}\label{subsec_byreflexivity} Let $T$ be any
locally fine and saturated fan, and $t\in T$ any point. By
reflexivity (proposition \ref{prop_reflex-dual}(iv)), the elements
$s\in\cO_{T,t}^\gp\otimes_\Z\Q$ correspond bijectively to
$\Q$-linear forms $\rho_s:U(t)(\Q_+)^\gp\to\Q$, and
$s\in\cO_{T,t}^\gp$ if and only if $\rho_s$ restricts to a morphism
of monoids $U(t)(\N)\to\Z$. Moreover, this bijection is compatible
with specialization maps : if $t'$ is a generization of $t$ in $T$,
then the form $U(t')(\Q_+)^\gp\to\Q$ induced by the image of $s$ in
$\cO_{T,t'}^\gp\otimes_\Z\Q$ is the restriction of $\rho_s$ (see
\eqref{subsec_from-face-to-prime}).

Hence, any global section $\lambda\in\Gamma(T,\cO_T^\gp)$ yields
a well defined function
$$
\rho_\lambda:T(\N)\to\Z
$$
whose restriction to $U(t)(\N)$ is the restriction of a $\Z$-linear
form on $U(t)(\N)^\gp$, for every $t\in T$; conversely, any such
function arises from a unique global section of $\cO_T^\gp$.
Likewise, we have a natural isomorphism between the $\Q$-vector
space of global sections $\lambda$ of $\cO_T^\gp\otimes_\Z\Q$, and
the space of functions $\rho_\lambda:T(\Q_+)\to\Q$ with a
corresponding linearity property.

Let now $\rho:T(\Q_+)\to\Q$ be any function; we may attach to $\rho$
a sheaf of fractional ideals of $\cO_{T,\Q}$ (notation of
\eqref{subsec_from-con-to-mon}), by the rule :
$$
\cI_{\rho,\Q}(U):=\{s\in\cO_{T,\Q}(U)~|~\rho_s\geq\rho_{|U}\} \qquad
\text{for every open subset $U\subset T$}
$$
In this generality, not much can be said concerning $\cI_{\rho,\Q}$;
to advance, we restrict our attention to a special class of functions,
singled out by the following :

\begin{definition}\label{def_roof}
Let $T$ be a locally fine and saturated fan.
\begin{enumerate}
\alphaenu
\item
A {\em roof\/} on $T$ is a function :
$$
\rho:T(\Q_+)\to\Q
$$
such that, for every $t\in T$, there exist $k:=k(t)\in\N$ and
$\Q$-linear forms
$$
\lambda_1,\dots,\lambda_k:U(t)(\Q_+)^\gp\to\Q
$$
with $\rho(s)=\min(\lambda_i(s)~|~i=1,\dots,k)$ for every
$s\in U(t)(\Q_+)$.
\item
An {\em integral roof\/} on $T$ is a roof $\rho$ on $T$ such that
$\rho(s)\in\Z$ for every $s\in T(\N)$.
\end{enumerate}
\end{definition}

\sset\subsubsection{}\label{subsec_encodes-geom} The interest of the
notion of roof on a fan $T$, is that it encodes in a geometrical
way, an integral subdivision of $T$, together with a coherent sheaf
of fractional ideals of $\cO_T$ (see definition
\ref{def_coh-idea-log}(iii)). This shall be seen in several steps.
To begin with, let $T$ and $\rho$ be as in definition
\ref{def_roof}(a). For any $t\in T$, pick a system
$\underline\lambda:=\{\lambda_1,\dots,\lambda_k\}$ of $\Q$-linear
forms fulfilling condition (b) of the definition; then for every
$i=1,\dots,k$ let us set :
$$
U(t,i)(\N):=\{x\in U(t)(\N)~|~\rho(x)=\lambda_i(x)\}.
$$
Notice that $U(t)(\N)=\cO_{T,t}$, by proposition
\ref{prop_reflex-dual}(iv). Moreover, say that $\underline\lambda$ is
{\em irredundant for $t$} if no proper subsystem of
$\underline\lambda$ fulfills condition (b) of definition
\ref{def_roof} relative to $U(t)$.

\begin{lemma}\label{lem_irredondo}
With the notation of \eqref{subsec_encodes-geom}, the following holds :
\begin{enumerate}
\item
$U(t,i)(\N)$ is a saturated fine monoid for every $i\leq k$.
\item
There is a unique system of\/ $\Q$-linear forms which is irredundant
for $t$.
\item
If $\underline\lambda$ is irredundant for $t$, then\ \
$\dim U(t,i)(\N)=\hgt_T(t)$\ \  for every $i=1,\dots,k$.
\end{enumerate}
\end{lemma}
\begin{proof}(i): We leave to the reader the verification that
$U(t,i)$ is a saturated monoid. Next, let
$\sigma_i\subset U(t)(\R_+)^{\gp\vee}$ be the convex polyhedral
cone spanned by the linear forms
$$
((\lambda_j-\lambda_i)\otimes_\Q\R~|~j=1,\dots,k).
$$
Then $\sigma_i$ is $\cO_{T,t}^{\gp\vee}$-rational, so $\sigma^\vee_i$
is $\cO_{T,t}^\gp$-rational, and $\sigma^\vee_i\cap\cO^\gp_{T,t}$
is a fine monoid (propositions \ref{prop_was-part-of-Gordon}(i),
and \ref{prop_Gordon}(i)), therefore the same holds for
$U(t,i)(\N)=\sigma^\vee_i\cap\cO_{T,t}$ (corollary
\ref{cor_fibres-are-fg}).

(iii): Notice that $\hgt_T(t)=\dim U(t)(\N)$, by proposition
\ref{prop_reflex-dual}(ii) and \eqref{eq_extisom-fans}. In view of
corollary \ref{cor_consequent}(i), it follows already that $\dim
U(t,i)\leq\hgt_T(t)$. Now, let $\lambda_1,\dots,\lambda_k$ be an
irredundant system, and suppose, by contradiction, that
$\dim\,U(t,i)(\N)<\dim\,U(t)(\N)$ for some $i\leq k$. Especially,
$\sigma^\vee_i\cap U(t)(\R_+)$ does not span the $\R$-vector space
$U(t)(\R_+)^\gp$, and therefore the dual $\sigma_i+U(t)(\R)^\vee$ is
not strictly convex (corollary \ref{cor_strongly}). After
relabeling, we may assume that $i=1$. Hence there exist
$a_j,b_j\in\R_+$ and $\phi,\phi'\in U(t)(\R)^\vee$, and an identity
:
$$
\sum_{j=2}^k a_j(\lambda_j-\lambda_i)+\phi= -\sum_{j=2}^k
b_j(\lambda_j-\lambda_i)-\phi'.
$$
Moreover, $\sum_{j=2}^k(a_j+b_j)>0$. It follows that there exist
$\psi\in U(t)(\R)^\vee$, and non-negative real numbers
$(c_j~|~j=2,\dots k)$ such that
$$
\lambda_i=\sum_{j=2}^k c_j\lambda_j+\psi \quad\text{and}\quad
\sum_{j=2}^kc_j=1.
$$
On the other hand, the irredundancy condition means that there
exists $x\in U(t,1)(\N)$ such that $\lambda_j(x)>\lambda_1(x)$ for
every $j>1$. Since $\psi(x)\geq 0$, we get a contradiction.

(ii): The assertion is clear, if $\hgt_T(t)\leq 1$. Hence suppose
that the height of $t$ is $\geq 2$, and let
$\underline\lambda:=\{\lambda_1,\dots,\lambda_k\}$ and
$\underline\mu:=\{\mu_1,\dots,\mu_r\}$ be two irredundant systems
for $t$. Fix $i\leq k$, and pick $x\in U(t,i)(\N)$ which does not
lie on any proper face of $U(t,i)(\N)$ (the existence of $x$ is
ensured by (iii) and proposition \ref{prop_Gordie}(i)); say that
$\mu_1(x)=\rho(x)$. Since $i$ is arbitrary, the assertion shall
follow, once we have shown that $\mu_1$ agrees with $\lambda_i$ on
the whole of $U(t,i)(\N)$.

However, by definition we have $\mu_1(y)\geq\lambda_i(y)$
for every $y\in U(t,i)(\N)$, and then it is easily seen that
$\Ker(\mu_1-\lambda_i)\cap U(t,i)(\N)$ is a face of $U(t,i)(\N)$;
since $x\in\Ker(\mu_1-\lambda_i)$, we deduce that $\mu_1-\lambda_i$
vanishes identically on $U(t,i)(\N)$.
\end{proof}

\sset\subsubsection{}\label{subsec_hencefan} Henceforth, we denote
by $\underline\lambda(t):=\{\lambda_1,\dots,\lambda_k\}$ the
irredundant system of $\Q$-linear forms for $t$. Let $1\leq i,j\leq
k$; then we claim that $U(t,i,j)(\N):=U(t,i)(\N)\cap U(t,j)(\N)$ is
a face of both $U(t,i)(\N)$ and $U(t,j)(\N)$. Indeed say that
$x,x'\in U(t,i)(\N)$ and $x+x'\in U(t,i,j)(\N)$; these conditions
translate the identities :
$$
\lambda_i(x)+\lambda_i(x')=\lambda_j(x)+\lambda_j(x') \qquad
\lambda_i(x)\leq\lambda_j(x) \qquad \lambda_i(x')\leq\lambda_j(x')
$$
whence $x,x'\in U(t,i,j)(\N)$. Define :
$$
U(t,i):=(\Spec\,U(t,i)(\N)^\vee)^\sharp \qquad
U(t,i,j):=(\Spec\,U(t,i,j)(\N)^\vee)^\sharp \qquad \text{for every
$i,j\leq k$}.
$$
According to \eqref{subsec_from-face-to-prime} and
\eqref{eq_extisom-fans}, the inclusion maps
$U(t,i,j)(\N)\to U(t,l)(\N)$ (for $l=i,j$) are dual to open
immersions
\set\begin{equation}\label{eq_immerse-roof}
U(t,i)\leftarrow U(t,i,j)\to U(t,j).
\end{equation}
We may then attach to $t$ and $\rho_{|U(t)}$ the fan $U(t,\rho)$
obtained by gluing the affine fans $U(t,i)$ along their common
intersections $U(t,i,j)$. The duals of the inclusions
$U(t,i)(\N)\to U(t)(\N)$ determine a well defined morphism of
locally fine and saturated fans :
\set\begin{equation}\label{eq_roof-rat-subdiv}
U(t,\rho)\to U(t)
\end{equation}
which, by construction, induces a bijection on $\N$-points :
$U(t,\rho)(\N)\isom U(t)(\N)$, so it is a rational subdivision,
according to proposition \ref{prop_val-crit-fans}.

\sset\subsubsection{}\label{subsec_global-roof}
Let now $t'\in T$ be a generization of $t$; clearly the system
$\underline\lambda':=\{\lambda'_1,\dots,\lambda'_k\}$ consisting
of the restrictions $\lambda'_i$ of the linear forms $\lambda_i$
to the $\Q$-vector subspace $U(t')(\Q_+)^\gp$, fulfills condition
(b) of definition \ref{def_roof}, relative to $t'$. However,
$\underline\lambda'$ may fail to be irredundant; after relabeling,
we may assume that the subsystem $\{\lambda'_1,\dots,\lambda'_l\}$
is irredundant for $t'$, for some $l\leq k$.
With the foregoing notation, we have obvious identities :
$$
U(t',i)(\N)=U(t,i)(\N)\cap U(t')(\N)
\qquad
U(t',i,j)(\N)=U(t,i,j)(\N)\cap U(t')(\N)
$$
for every $i,j\leq l$; whence, in light of remark
\ref{rem_reflex-dual}(iii), a commutative diagram of fans :
$$
\xymatrix{
U(t',i) \ar[d] & \ar[l] U(t',i,j) \ar[r] \ar[d] & U(t',j) \ar[d] \\
U(t,i)\times_{U(t)}U(t') &
\ar[l] U(t,i,j)\times_{U(t)}U(t') \ar[r] &
U(t,i)\times_{U(t)}U(t') }
$$
whose top horizontal arrows are the open immersions
\eqref{eq_immerse-roof} (with $t$ replaced by $t'$),
whose bottom horizontal arrows are the open immersions
$\eqref{eq_immerse-roof}\times_{U(t)}U(t')$, and whose vertical
arrows are natural isomorphisms. Since $U(t')$ is an open subset
of $U(t)$, we deduce an open immersion
$$
j_{t,t'}:U(t',\rho)\to U(t,\rho).
$$
If $t''$ is a generization of $t'$, it is clear that
$j_{t',t''}\circ j_{t,t'}=j_{t,t''}$, hence we may glue the
fans $U(t,\rho)$ along these open immersions, to obtain a
locally fine and saturated fan $T(\rho)$. Furthermore, the
morphisms \eqref{eq_roof-rat-subdiv} glue to a single rational
subdivision :
\set\begin{equation}\label{eq_single}
T(\rho)\to T.
\end{equation}

\begin{remark}\label{rem_roofs-dominate}
(i)\ \ 
In the language of definition \ref{def_fans} the foregoing
lengthy procedure translates as the following simple geometric
operation. Given a fan $\Delta$ (consisting of a collection of
convex polyhedral cones of a $\R$-vector space $V$), a roof on
$\Delta$ is a piecewise linear function
$F:=\bigcup_{\sigma\in\Delta}\sigma\to\R$, which is concave on each
$\sigma\in\Delta$ (and hence it is a roof on each such $\sigma$,
in the sense of example \ref{ex_roof-again}).
Then, such a roof determines a natural refinement $\Delta'$ of
$\Delta$; namely, $\Delta'$ is the coarsest refinement such that,
for each $\sigma'\in\Delta'$, the function $\rho_{|\sigma'}$ is
the restriction of a $\R$-linear form on $V$.
This refinement $\Delta'$ corresponds to the present $T(\rho)$.

(ii)\ \ 
Moreover, let $P$ be a fine, sharp and saturated monoid of
dimension $d$, set $T_P:=(\Spec\,P)^\sharp$, and suppose that
$f:T\to T_P$ is any integral, fine, proper and saturated
subdivision. Then $f$ corresponds to a geometrical subdivision
$\Delta$ of the strictly convex polyhedral cone
$T_P(\R_+)=P^\vee_\R$, and we claim that $\Delta$ can be refined by
the subdivision associated with a roof on $T_P$. Namely, let
$\Delta_{d-1}$ be the subset of $\Delta$ consisting of all $\sigma$
of dimension $d-1$; every $\sigma\in\Delta_{d-1}$ is the
intersection of a $d$-dimensional element of $\Delta$ and a
hyperplane $H_\sigma\subset P^{\gp\vee}_\R$; such hyperplane is the
kernel of a linear form $\lambda_\sigma$ on $P^{\gp\vee}_\R$. Let us
define
$$
\rho(x):=\sum_{\sigma\in\Delta_{d-1}}\min(0,\lambda_\sigma(x))
\qquad \text{for every $x\in T_P(\R_+)$.}
$$
Then it is easily seen that the subdivision of $T_P(\R_+)$
associated with the roof $\rho$ as in (i), refines the subdivision
$\Delta$. In the language of fans, this construction translates as
follows. For every point $\sigma\in T$ of height $d-1$, let
$H_\sigma\subset P^\gp$ be the kernel of the surjection
$P^\gp\to\cO^\gp_{T,\sigma}$ induced by $\log f$; notice that
$H_\sigma$ is a free abelian group of rank one, and pick a generator
$s_\sigma$ of $H_\sigma$, which -- as in
\eqref{subsec_byreflexivity} -- corresponds to a function
$\lambda_\sigma:T_P(\N)\to\Z$, so we may again consider the integral
roof $\rho$ on $T_P$ defined as in the foregoing. Then it is easily
seen that the rational subdivision $T(\rho)\to T_P$ associated with
$\rho$, factors as the composition of $f$ and a (necessarily unique)
integral subdivision $g:T(\rho)\to T$. More precisely, for every
mapping $\eps:\{\sigma\in T~|~\hgt(\sigma)=d-1\}\to\{0,1\}$, let us
set
$$
\lambda_\eps:=\sum_{\hgt(\sigma)=d-1}\eps(\sigma)\cdot\lambda_\sigma
\qquad\text{and}\qquad U(\eps)(\N):=\{x\in
T_P(\N)~|~\rho(x)=\lambda_\eps(x)\}.
$$
Whenever $U(\eps)(\N)$ has dimension $d$, let $t_\eps\in T(\rho)$ be
the unique point such that $U(\eps)(\N)^\vee=\cO_{T(\rho),t_\eps}$.
As the reader may check, there exists a unique closed point $\tau\in
T$, such that $\eps(\sigma)\cdot\lambda_\sigma(x)\leq 0$ for every
$x\in U(\tau)(\N)$ and every $\sigma\in U(\tau)$ of height $d-1$.
Then we have $g(t_\eps)=\tau$, and the restriction $U(t_\eps)\to
U(\tau)$ of $g$ is deduced from the inclusion $U(\eps)(\N)\subset
U(\tau)(\N)$ of submonoids of $T_P(\N)$.

(iii)\ \ Furthermore, in the situation of (ii), the roof $\rho$ on
$T_P$ can also be viewed as a roof on $T$, and then it is clear from
the construction that the morphism $g:T(\rho)\to T$ is also the
subdivision of $T$ attached to the roof $\rho$.
\end{remark}

We wish now to establish some basic properties of the sheaf
of fractional ideals
$$
\cI_{\!\rho}:=\cI_{\!\rho,\Q}\cap\cO_T^\gp
$$
attached to a given roof on $T$. First we remark :

\begin{lemma}\label{lem_coh-roof}
Keep the notation of \eqref{subsec_hencefan}, and let
$s\in\cO^\gp_{T,t}\otimes_\Z\Q$ be any element such that
$\rho_s\geq\rho_{|U(t)}$. Then we have :
\begin{enumerate}
\item
There exist $\phi\in(\cO_{T,t})_\Q$ (notation of
\eqref{subsec_from-con-to-mon}) and $c_1,\dots,c_k\in\R_+$ such that
:
$$
\rho_s=\rho_\phi+\sum_{i=1}^k c_i\lambda_i \qquad \sum_{i=1}^kc_i=1.
$$
\item
The stalk $\cI_{\!\rho,t}$ is a finitely generated $\cO_{T,t}$-module.
\end{enumerate}
\end{lemma}
\begin{proof} (i): Let $\sigma\subset U(t)(\R_+)^\gp$ be the
convex polyhedral cone spanned by the linear forms
$((\lambda_i-\rho_s)\otimes_\Q\R~|~i=1,\dots,k)$. Then the assumption
on $s$ means that
$$
U(t)(\R_+)\cap\sigma^\vee=U(t)(\R_+)\cap\Ker\,\rho_s\otimes_\Q\R.
$$
Especially $\sigma^\vee\cap U(t)(\R_+)$ does not span $U(t)(\R_+)^\gp$.
Then one can repeat the proof of lemma \ref{lem_irredondo}(iii)
to derive the assertion.

(ii): By remark \ref{rem_reflex-dual}(i), we may write
$\lambda_i=\rho_{s_i}-\rho_{s'_i}$, where
$s_i,s'_i\in(\cO_{T,t})_\Q$ for each $i\leq k$; pick $N\in\N$ large
enough, so that $Ns'_i\in\cO_{T,t}$ for every $i\leq k$, and set
$\tau:=N\sum_{i=1}^ks'_i$. Then
$\tau+\cI_{\!\rho,t}\subset\cO^\gp_{T,t}\cap(\cO_{T,t})_\Q=\cO_{T,t}$.
By proposition \ref{prop_ideals-in-fg-mon}(ii), we deduce that
$\tau+\cI_{\!\rho,t}$ is a finitely generated ideal, whence the
contention.
\end{proof}

\begin{proposition} Let $T$ be a locally fine and saturated
fan, $\rho$ a roof on $T$. Then the associated fractional
ideal $\cI_{\!\rho}$ of $\cO_T$ is coherent.
\end{proposition}
\begin{proof} In view of lemma \ref{lem_coh-roof}(ii) and remark
\ref{rem_local-trivial-qcoh}(v), it suffices to show that
$\cI_{\!\rho}$ is quasi-coherent, {\em i.e.} for every
generization $t'$ of $t$, the image of
$\cI_{\!\rho,t}$ in $\cO_{T,t'}$ generates the $\cO_{T,t'}$-module
$\cI_{\!\rho,t'}$. Fix such $t'$; by propositions
\ref{prop_was-part-of-Gordon}(i) and \ref{prop_Gordie}(i), there
exists $\lambda\in\cO_{T,t}=U(t)(\N)^\vee$ such that
$U(t')(\N)=\Ker\,\lambda$; especially, we see that $U(t')(\N)^\gp$
is a direct summand of $U(t)(\N)^\gp$. Now, let
$s'\in\cI_{\!\rho,t'}$ be any local section; it follows that we may
find $s\in\cO^\gp_{T,t}$ such that $\rho_s:U(t)(\N)^\gp\to\Z$ is a
$\Z$-linear extension of the corresponding $\Z$-linear form
$\rho_{s'}:U(t')(\N)^\gp\to\Z$. Let also
$\{\lambda_1,\dots,\lambda_k\}$ be the irredundant system of
$\Q$-linear forms for $t$ (relative to the roof $\rho$). For every
$i\leq k$ we have the following situation :
$$
\lambda_\R:=\lambda\otimes_\Q\R\in U(t,i)(\R_+)^\vee \qquad
(\rho_s-\lambda_i)\otimes_\Q\R\in U(t',i)(\R_+)^\vee.
$$
However, $U(t',i)(\R_+)=U(t,i)(\R_+)\cap\Ker\,\lambda_\R$, hence
$U(t',i)(\R_+)^\vee=U(t,i)(\R_+)^\vee+\R\lambda_\R$. Especially,
there exist $r_i\in\R_+$ and $\phi\in U(t)(\R_+)^\vee$ such that
$(\rho_s-\lambda_i)\otimes_\Q\R=\phi-r_i\lambda_\R$. Let $N$ be an
integer greater than $\max(r_1,\dots,r_k)$; it follows that
$s+N\lambda\in\cI_{\!\rho,t}$ and its image in $\cI_{\!\rho,t'}$
equals $s'$.
\end{proof}

\begin{proposition}\label{prop_to-subdiv-is-to-blow}
In the situation of \eqref{subsec_encodes-geom} suppose that
$\rho$ is an integral roof on $T$. Then the morphism
\eqref{eq_single} is the saturation of a blow up of the
fractional ideal $\cI_{\!\rho}$.
\end{proposition}
\begin{proof} We have to exhibit an isomorphism
$f:T(\rho)\to X:=\Proj\,\cB(\cI)^\sat$ of $T$-fans. We begin with :

\begin{claim}\label{cl_identify-sat} (i)\ \ The fractional ideal
$\cI_{\!\rho}\cO_{T(\rho)}$ is invertible.
\begin{enumerate}
\addenu
\item
For every $n\in\N$, denote by $n\rho:T(\Q_+)\to\Q_+$ the function
given by the rule $x\mapsto n\cdot\rho(x)$ for every $x\in T(\Q_+)$.
Then :
$$
\cB(\cI)^\sat=\cB':=\coprod_{n\in\N}\cI_{\!n\rho}
$$
\end{enumerate}
\end{claim}
\begin{pfclaim}(ii): The assertion is local on $T$, hence we may
assume that $T=U(t)$ for some $t\in T$, in which case, denote by
$\underline\lambda:=\{\lambda_1,\dots,\lambda_k\}$ the irredundant
system of $\Q$-linear forms for $t$. Let $n\in\N$, and
$s\in\cI_{n\rho}(U(t))$; it is easily seen that the $T$-monoid
$\cB'$ is saturated, hence it suffices to show that there exists an
integer $k>0$ such that $k\rho_s=\rho_{s'}$ for some
$s'\in\cI^k(U(t))$ (notation of \eqref{subsec_byreflexivity}).
However, lemma \ref{lem_coh-roof}(ii) implies more precisely
that we may find such $k$, so that the corresponding $s'$
lies in the ideal generated by $\underline\lambda$.

(i): The assertion is local on $T(\rho)$, hence we consider
again $t\in T$ and the corresponding $\underline\lambda$ as in
the foregoing. It suffices to show that $\cJ:=\cI_{\!\rho}\cO_{U(t,i)}$
is invertible for every $i=1,\dots,k$ (notation of
\eqref{subsec_hencefan}). However, by inspecting the constructions
it is easily seen that $\cJ(U(t,i))$ consists of all
$s\in(U(t,i)(\N)^\vee)^\gp$ such that $\rho_s(x)\geq\rho(x)$
for every $x\in U(t,i)(\Q_+)$, {\em i.e.} $\rho_s(x)\geq\lambda_i(x)$
for every $x\in U(t,i)(\Q_+)$. However, since $\rho$ is integral,
we have $\lambda_i\in\cO_{T,t}^\gp$; if we apply lemma
\ref{lem_coh-roof}(i) with $T$ replaced by $U(t,i)$, we conclude that
$\cJ(U(t,i))$ is the fractional ideal generated by $\lambda_i$,
whence the contention.
\end{pfclaim}

In view of claim \ref{cl_identify-sat}(i) we see that there
exists a unique morphism $f$ of $T$-fans from $T(\rho)$ to $X$. It
remains to check that $f$ is an isomorphism; the latter assertion
is local on $X$, hence we may assume that $T=U(t)$ for some
$t\in T$, and then we let again $\underline\lambda$ be the
irredundant system for $t$. A direct inspection yields a natural
identification
of $\cO_{T,t}$-monoids :
$$
\cB'(U(t,i))_{(\lambda_i)}\isom U(t,i)(\N)^\vee
\qquad\text{for every $i=1,\dots,k$}
$$
whence an isomorphism $U(t,i)\isom D_+(\lambda_i)\subset X$,
which -- by uniqueness -- must coincide with the restriction
of $f$. On the other hand, the proof of claim
\ref{cl_identify-sat}(ii) also shows that
$X=D_+(\lambda_1)\cup\cdots\cup D_+(\lambda_k)$, and the
proposition follows.
\end{proof}

\begin{example}\label{ex_cut-by-hyperplane}
Let $P$ be a fine, sharp and saturated monoid.

(i)\ \ The simplest non-trivial roofs on $T_P:=(\Spec\,P)^\sharp$
are the functions $\rho_\lambda$ such that
$$
\rho_\lambda(x):=\min(0,\lambda(x)) \qquad\text{for every $x\in
T_P(\Q_+)$.}
$$
where $\lambda$ is a given element of
$\Hom_\Q(T_P(\Q_+)^\gp,\Q)\simeq P^\gp\otimes_\Z\Q$. Such a
$\rho_\lambda$ is an integral roof, provided $\lambda\in P^\gp$. In
the latter case, we may write $\lambda=\rho_{s_1}-\rho_{s_2}$, for
some $s_1,s_2\in P$. Set $\rho':=\rho_\lambda+\rho_{s_2}$, {\em
i.e.} $\rho'=\min(\rho_{s_1},\rho_{s_2})$mplest non-trivial roofs on $T_P:=(\Spec\,P)^\sharp$
are the functions $\rho_\lambda$ such that
$$
\rho_\lambda(x):=\min(0,\lambda(x)) \qquad\text{for every $x\in
T_P(\Q_+)$.}
$$
where $\lambda$ is a given element of
$\Hom_\Q(T_P(\Q_+)^\gp,\Q)\simeq P^\gp\otimes_\Z\Q$. Such a
$\rho_\lambda$ is an integral roof, provided $\lambda\in P^\gp$. In
the latter case, we may write $\lambda=\rho_{s_1}-\rho_{s_2}$, for
some $s_1,s_2\in P$. Set $\rho':=\rho_\lambda+\rho_{s_2}$, {\em
i.e.} $\rho'=\min(\rho_{s_1},\rho_{s_2})$; clearly
$T(\rho_\lambda)=T(\rho')$, and on the other hand lemma
\ref{lem_coh-roof}(i) implies that the ideal $\cI_{\rho'}$ is the
saturation of the ideal generated by $s_1$ and $s_2$.

(ii)\ \ More generally, any system
$\lambda_1,\dots,\lambda_n\in\Hom_\Q(T_P(\Q_+)^\gp,\Q)$ of
$\Q$-linear forms yields a roof $\rho$ on $T_P$, such that
$\rho(x):=\sum_{i=1}^n\min(0,\lambda_i(x))$ for every $x\in
T_P(\Q_+)$. A simple inspection shows that the corresponding
subdivision $T(\rho)\to T$ can be factored as the composition of $n$
subdivisions $g_i:T_i\to T_{i-1}$, where $T_0:=T_P$, $T_n:=T(\rho)$,
and each $g_i$ (for $i\leq n$) is the subdivision of $T_{i-1}$
corresponding to the roof $\rho_{\lambda_i}$ as defined in (i).

(iii)\ \ These subdivisions of $T_P$ "by hyperplanes" are precisely
the ones that occur in remark \ref{rem_roofs-dominate}(ii),(iii).
Summing up, we conclude that every proper integral and saturated
subdivision $g:T\to T_P$ of $T_P$ can be dominated by another
subdivision $f:T(\rho)\to T_P$ of the type considered in (ii), so
that $f$ factors as the composition of $g$ and a subdivision
$h:T(\rho)\to T$ which is also of the type (ii). Especially, both
$f$ and $h$ can be realized as the composition of finitely many
saturated blow up of ideals generated by at most two elements of
$P$.
\end{example}

\sset\subsubsection{}\label{subsec_std-subdivision}
Let $P$ be a fine, sharp and saturated monoid. A proper,
integral, fine and saturated subdivision of
$$
T_P:=(\Spec\,P)^\sharp
$$
is essentially equivalent to a $(P^\gp)^\vee$-rational
subdivision of the polyhedral cone $\sigma:=P_\R^\vee$
(see \eqref{subsec_Gordon-more} and definition \ref{def_fans}).
A standard way to subdivide a polyhedron $\sigma$ consists in
choosing a point $x_0\in\sigma\setminus\!\{0\}$, and forming all the
polyhedra $x_0*F$, where $F$ is any proper face of $\sigma$, and
$x_0*F$ denotes the convex span of $x_0$ and $F$.
We wish to describe the same operation in terms of the topological
language of affine fans.

Namely, pick any non-zero $\phi\in T_P(\Q_+)$
($\phi$ corresponds to the point $x_0$ in the foregoing). Let
$U(\phi)\subset\Spec\,P$ be the set of all prime ideals $\fp$ such
that $\phi(P\!\setminus\!\fp)\neq\{0\}$; in other words, the complement
of $U(\phi)$ is the topological closure of the support of $\phi$
in $T_P$, especially, $U(\phi)$ is an open subset of $T_P$.
Denote by $j:P=\Gamma(T_P,\cO_{T_P})\to\Gamma(U(\phi),\cO_{T_P})$ the
restriction map. The morphism of monoids :
$$
P\to\Gamma(U(\phi),\cO_{T_P})\times\Q_+ \qquad
x\mapsto(j(x),\phi(x))
$$
determines a cocartesian diagram of fans
$$
\xymatrix{ U(\phi)\times(\Spec\,\Q_+)^\sharp \ar[d]_{\psi'}
\ar[rr]^-{\beta'} & & U(\phi)\times(\Spec\,\N)^\sharp \ar[d]^\psi
\\
T_P \ar[rr]^-\beta & & T_{\phi^{-1}\N}:=(\Spec\,\phi^{-1}\N)^\sharp
}
$$

\begin{lemma}\label{lem_star-subdivide}
With the notation of \eqref{subsec_std-subdivision}, the morphisms
$\psi$ and $\psi'$ are proper rational subdivisions, which we call
the {\em subdivisions centered at $\phi$}.
\end{lemma}
\begin{proof} Notice that both $\beta$ and $\beta'$ are
homeomorphisms on the underlying topological spaces, and moreover
both $\log\beta^\gp\otimes_\Z\one_\Q$ and
$\log\beta^{\prime\gp}\otimes_\Z\one_\Q$ are isomorphisms. Thus, it
suffices to show that $\psi$ is a proper rational subdivision, hence
we may replace $P$ by $\phi^{-1}\N$, which allows to assume that
$\phi$ is a morphism of monoids $P\to\N$, and $\psi$ is a morphism
of fans
$$
T':=U(\phi)\times(\Spec\,\N)^\sharp\to T_P.
$$
In this situation, by inspecting the construction, we find that
$\psi$ restricts to an isomorphism :
$$
\psi^{-1}U(\phi)\isom U(\phi)
$$
and the preimage of the closed subset $T_P\!\setminus\!U(\phi)$ is
the preimage of the closed point $\fm\in(\Spec\,\N)^\sharp$ under
the natural projection $T'\to(\Spec\,\N)^\sharp$; in view of the
discussion of \eqref{prod_of-fine-fans}, this is naturally
identified with $U(\phi)\times\{\fm\}$. Moreover, the restriction
$U(\phi)\times\{\fm\}\to T_P\setminus U(\phi)$ of $\psi$ is the map
$(t,\fm)\mapsto t\cup\phi^{-1}\fm$ (recall that $t\subset P$ is a
prime ideal which does not contain $\phi^{-1}\fm$). Thus, the
assertion will follow from :
\begin{claim} The $\Q$-linear map
$$
\log\psi^\gp_{(t,\fm)}\otimes_\Z\one_\Q:
P^\gp\otimes_\Z\Q\to(\cO^\gp_{T_P,t}\times\Z)\otimes_\Z\Q
$$
is surjective for every $t\in U(\phi)$, and the induced map :
\set\begin{equation}\label{eq_induced-on-morphs}
T'(\Q_+)\to T_{\phi^{-1}\N}(\Q_+)=T_P(\Q_+)
\end{equation}
is bijective.
\end{claim}
\begin{pfclaim}[] Indeed, let $F_t\subset P^\vee$ be the face of $P^\vee$
corresponding to the point $t\in U(\phi)$, under the bijection
\eqref{eq_face-to-face}; then $\phi\notin F_t$, whence a natural
isomorphism of monoids :
$$
(F_t+\N\phi)^\vee\isom\cO_{T_P,t}\times\N \quad : \quad
\lambda\mapsto(\lambda_{|F},\lambda(\phi))
$$
whose inverse, composed with $\log\psi_{(t,\fm)}$, yields the
restriction map $P\to(F_t+\N\phi)^\vee$. This interpretation makes
evident the surjectivity of
$\log\psi^\gp_{(t,\fm)}\otimes_\Z\one_\Q$. In view of example
\ref{ex_fans}(iii), the bijectivity of \eqref{eq_induced-on-morphs}
is also clear, if one remarks that :
$$
P^\vee_\R=\bigcup_{t\in U(\phi)}(F_t+\N\phi)_\R.
$$
The latter identity is obvious from the geometric interpretation in
terms of polyhedral cones. A formal argument runs as follows. Let
$\phi'\in P^\vee_\R$; since $P$ spans $P^\gp_\R$, the cone
$(P_\R)^\vee$ is strictly convex (corollary \ref{cor_strongly}),
hence the line $\phi'+\R\phi\subset P^\gp_\R$ is not contained in
$P^\vee_\R$, therefore there exists a largest $r\in\R$ such that
$\phi'-r\phi\in P^\vee_\R$, and necessarily $r\geq 0$. If
$\phi'-r\phi=0$, the assertion is clear; otherwise, let $F$ be the
minimal face of $P^\vee$ such that $\phi'-r\phi\in F_\R$, so that
$\phi'=(\phi'-r\phi)+r\phi\in(F+\N\phi)_\R$. Thus, we are reduced to
showing that $\phi\notin F$. But notice that $\phi'-r\phi$ lies in
the relative interior of $F$; therefore, if $\phi\in F$, we may find
$\eps>0$ such that $\phi'-(r+\eps)\phi$ still lies in $F_\R$,
contradicting the definition of $r$.
\end{pfclaim}
\end{proof}

\sset\subsubsection{}\label{subsec_extend-roofs}
Lemma \ref{lem_star-subdivide} is frequently used to construct
subdivisions centered at an {\em interior point\/} of $T_P$,
{\em i.e.} a point $\phi\in T_P(\N)$ which does not lie on any
proper face of $T_P(\N)$ (equivalently, the support of $\phi$ is
the closed point $\fm_P$ of $T_P$). In this case
$U(\phi)=T_P\!\setminus\!\{\fm_P\}=(T_P)_{d-1}$, where $d:=\dim\,P$.
By lemma \ref{lem_star-subdivide}, the $\N$-point $\phi$ lifts
to a unique $\Q_+$-point $\tilde\phi$ of
$(T_P)_{d-1}\times(\Spec\,\N)^\sharp$, and by inspecting
the definitions, it is easily seen that -- under the identification
of remark \ref{rem_product-of-aff-fans}(iii) -- the support of
$\tilde\phi$ is the point $(\emptyset,\fm_\N)$, where
$\emptyset\in T_P$ is the generic point. More precisely, we may
identify $(\Spec\,\N)^\sharp(\Q_+)$ with $\Q_+$, and
$(T_P)_{d-1}(\Q_+)$ with a cone in the $\Q$-vector space
$T_P(\Q_+)^\gp$, and then $\tilde\phi$ corresponds to the
point $(0,1)\in T_P(\Q_+)^\gp\times\Q_+$.

Suppose now that we have an integral roof
$$
\rho:(T_P)_{d-1}(\Q_+)\to\Q
$$
and let $\pi:T(\rho)\to(T_P)_{d-1}$ be the associated subdivision.
Let also $\psi:(T_P)_{d-1}\times(\Spec\,\N)^\sharp\to T_P$ be
the subdivision centered at $\phi$; we deduce a new subdivision :
\set\begin{equation}\label{eq_combine-subdivisions}
T^*:=T(\rho)\times(\Spec\,\N)^\sharp
\xrightarrow{ \pi\times\one_{(\Spec\,\N)^\sharp} }
(T_P)_{d-1}\times(\Spec\,\N)^\sharp\xrightarrow{ \psi } T_P
\end{equation}
whose restriction to the preimage of $(T_P)_{d-1}$ is $T_P$-isomorphic
to $\pi$.

\begin{lemma}\label{lem_extend-roof}
In the situation of \eqref{subsec_extend-roofs}, there exist an
integral roof
$$
\tilde\rho:T_P(\Q_+)\to\Q
$$
whose restriction to $(T_P)_{d-1}(\Q_+)$ agrees with $\rho$, and
a morphism $T^*\to T(\tilde\rho)$ of $T_P$-monoids, whose underlying
continuous map is a homeomorphism.
\end{lemma}
\begin{proof} According to remark \ref{rem_product-of-aff-fans}(vi),
we have a natural identification :
$$
T_P(\Q_+)=T(\rho)(\Q_+)\times(\Spec\,\N)^\sharp(\Q_+)=
T(\rho)(\Q_+)\times\Q_+\subset T_P(\Q_+)^\gp\times\Q.
$$
mapping the point $\phi$ to $(0,1)$. For a given $c\in\R$,
denote by $\rho_c:T_P(\Q_+)\to\Q$ the function given by the
rule : $(x,y)\mapsto\rho(x)+cy$ for every $x\in T(\rho)(\Q_+)$
and every $y\in\Q_+$.

Let $t\in T(\rho)$ be any point of height $d-1$; by assumption,
there exists a $\Q$-linear form $\lambda_t:U(\pi(t))(\Q_+)^\gp\to\Q$
whose restriction to $U(t)(\Q_+)$ agrees with the restriction of
$\rho$. Therefore, the restriction of $\rho_c$ to
$U(t,\fm_\N)(\Q_+)=U(t)(\Q_+)\times\Q_+$ agrees with the
restriction of the $\Q$-linear form
$$
\lambda_{t,c}:U(\pi(t))(\Q_+)^\gp\times\Q=T_P(\Q_+)^\gp\to\Q
\qquad
(x,y)\mapsto\lambda_t(x)+cy.
$$
For any two points $t,t'\in T(\rho)$ of height $d-1$, with
$\pi(t)=\pi(t')$, pick a finite system of generators
$\{x_1,\dots,x_n\}$ for $U(t')(\N)$, and let
$y_1,\dots,y_n\in U(t)(\Q)^\gp$, $a_1,\dots,a_n\in\Q$ such that
$$
x_i=y_i+a_i\phi
\qquad
\text{in the vector space $U(\pi(t))(\Q_+)^\gp$}.
$$
In case $x_i\in U(t)(\N)$, we shall have $a_i=0$, and otherwise
we remark that $a_i>0$. Indeed, if $a_i<0$ we would have
$x_i-a_i\phi\in U(t)(\Q_+)^\gp\cap T_P(\Q_+)\subset U(\pi(t))(\Q_+)$;
however, $U(\pi(t))(\Q_+)$ is a proper face of $T_P(\Q_+)$, hence
$\phi\in U(\pi(t))(\Q_+)$, a contradiction.

We may then find $c>0$ large enough, so that
$\lambda_{t'}(x_i)<\lambda_{t,c}(x_i)$ for every $x_i\notin U(t)(\N)$.
It follows easily that
\set\begin{equation}\label{eq_ineq-lambdas}
\lambda_{t',c}(x)<\lambda_{t,c}(x)
\qquad
\text{for all $x\in U(t',\fm_\N)(\Q_+)\!\setminus\!U(t,\fm_\N)(\Q_+)$}.
\end{equation}
Clearly we may choose $c$ large enough, so that \eqref{eq_ineq-lambdas}
holds for every pair $t,t'$ as above, and then it is clear that
$\rho_c$ will be a roof on $T_P$. Notice that the points
of height $d$ of $T^*$ are precisely those of the form $(t,\fm_\N)$,
for $t\in T_P$ of height $d-1$, so the points of $T(\rho_c)$ of
height $d$ are in natural bijection with those of $T^*$, and if
$\tau\in T(\rho_c)$ corresponds to $\tau^*\in T^*$ under this
bijection, we have an injective map $U(\tau^*)(\N)\to U(\tau)(\N)$,
commuting with the induced projections to $T_P(\N)$.
There follows a morphism of $T_P$-fans $U(\tau^*)\to U(\tau)$
inducing a bijection $U(\tau^*)(\Q_+)\isom U(\tau)(\Q_+)$.
Since $T(\rho_c)$ (resp. $T^*$) is the union of all such
$U(\tau)$ (resp. $U(\tau^*)$), we deduce a morphism
$T^*\to T(\rho_c)$ inducing a homeomorphism on underlying
topological spaces. Lastly, say that $\tau^*=(t,\fm)$; then
$U(\tau^*)(\N)^\gp=U(t)(\N)^\gp\oplus\Z\phi$, and on the other
hand $U(t)(\N)^\gp$ is a direct factor of $U(\tau)(\N)^\gp$
(since the specialization map
$\cO^\gp_{T(\rho_c),\tau}\to\cO^\gp_{T,t}$ is surjective).
It follows easily that we may choose for $c$ a suitable positive
integer, in such a way that the resulting roof $\tilde\rho:=\rho_c$
will also be integral.
\end{proof}

\sset\subsubsection{}
Let $P$ be a fine, sharp and saturated monoid, and set as usual
$T_P:=(\Spec\,P)^\sharp$. There is a canonical choice of a point
in $T_P(\N)$ which does not lie on any proper face of $T_P(\N)$;
namely, one may take the $\N$-point $\phi_P$ defined as the
sum of the generators of the one-dimensional faces of $T_P(\N)$
(such faces are isomorphic to $\N$, by theorem
\ref{th_structure-of-satu}(ii)).

Set $d:=\dim\,P$; if $\rho$ is a given integral roof for $(T_P)_{d-1}$,
and $c\in\N$ is a sufficiently large, we may then attach to
the datum $(\phi_P,\rho,c)$ an integral roof $\tilde\rho$ of $T_P$
extending $\rho$ as in lemma \ref{lem_extend-roof}, and such
that $\tilde\rho(\phi_P)=c$. More generally, let $T$ be a
fine and saturated fan of dimension $d$, and suppose we have
a given integral roof $\rho_{d-1}$ on $T_{d-1}$; for every point
$t\in T$ of height $d$ we have the corresponding canonical point
$\phi_t$ in the ``interior'' of $U(t)(\N)$, and we may then pick
an integer $c_d$ large enough, so that $\rho_{d-1}$ extends to an
integral roof $\rho_d:T(\Q_+)\to\Q$ with $\rho_d(\phi_t)=c_d$ for
every $t\in T$ of height $d$, and such that the associated
subdivision is $T_P$-homeomorphic to $T(\rho_{d-1})\times(\Spec\,\N)$.

This is the basis for the inductive construction of an integral
roof on $T_P$ which is canonical in a certain restricted sense.
Indeed, fix an increasing sequence of positive integers
$\underline c:=(c_2,\dots,c_d)$; first we define
$\rho_1:T_1(\Q_+)\to\Q_+$ to be the identically zero map.
This roof is extended recursively to a function $\rho_h$ on $T_h$,
for each $h=2,\dots,d$, by the rule given above, in such a way that
$\rho_h(\phi_t)=c_i$ for every point $t$ of height $i\leq h$.
By the foregoing we see that the sequence $\underline c$ can
be chosen so that $\rho_d$ shall again be an integral roof.

\sset\subsubsection{}\label{subsec_models-of-fans}
More generally, let $\cS:=\{P_1,\dots,P_k\}$ be any finite
set of fine, sharp and saturated monoids. We let $\cS\text{-}\Fan$
be the full subcategory of $\Fan$ whose objects are the fans
$T$ such that, for every $t\in T$, there exist $P\in\cS$ and
an open immersion $U(t)\subset(\Spec\,P)^\sharp$. Then the
foregoing shows that we may find a sequence of integers
$\underline c(\cS):=(c_2,\dots,c_d)$, with
$d:=\max(\dim\,P_i~|~i=1,\dots,k)$ such that the following holds.
Every object $T$ of $\cS\text{-}\Fan$ is endowed with an integral
roof $\rho_T:T(\Q_+)\to\Q_+$ such that :
\begin{itemize}
\item
$\rho_T(\phi_t)=c_i$ whenever $\hgt(t)=i\geq 2$, and $\rho_T$
vanishes on $T_1(\Q_+)$.
\item
Every open immersion $g:T'\to T$ in $\cS\text{-}\Fan$ determines
an open immersion $\tilde g:T'(\rho_{T'})\to T(\rho_T)$ such
that the diagram
$$
\xymatrix{
        T'(\rho_{T'}) \ar[r]^-{\tilde g} \ar[d]_{\pi_{T'}} &
        T(\rho_T) \ar[d]^{\pi_T} \\
        T' \ar[r]^-g & T
}$$
commutes (where $\pi_T$ and $\pi_{T'}$ are the subdivisions associated
to $\rho_T$ and $\rho_{T'}$).
\item
If $\dim\,T=d$, there exists a natural rational subdivision
\set\begin{equation}\label{eq_extra-subdiv}
T_{d-1}(\rho_{T_{d-1}})\times(\Spec\,\N)^\sharp\to T(\rho_T)
\end{equation}
which is a homeomorphism on the underlying topological spaces.
\end{itemize}

\sset\subsubsection{}\label{subsec_canonical-coordinates}
Notice as well that, by construction, $T_1(\rho_{T_1})=T_1$;
hence, by composing the morphisms \eqref{eq_extra-subdiv}, we
obtain a rational subdivision of $T(\rho_T)$ :
$$
T_1\times(\Spec\,\N^{\oplus d-1})^\sharp\to\cdots\to
T_{d-2}(\rho_{T_{d-2}})\times(\Spec\,\N^{\oplus 2})^\sharp\to
T_{d-1}(\rho_{T_{d-1}})\times(\Spec\,\N)^\sharp\to T(\rho_T).
$$
In view of theorem \ref{th_structure-of-satu}(ii), it is
easily seen that every affine open subsect of $T_1$ is
isomorphic to $(\Spec\,\N)^\sharp$, and any two such open subsets
have either empty intersection, or else intersect in their
generic points. In any case, we deduce a natural epimorphism
$$
(\Q_{+})_{T_1}\to\cO_{T_1,\Q}
$$
(notation of \eqref{subsec_from-con-to-mon}) from the constant
$T_1$-monoid arising from $\Q_+$; whence an epimorphism :
$$
\theta_T:(\Q^{\oplus d}_+)_{T(\rho_T)}\to\cO_{T(\rho_T),\Q}
$$
which is compatible with open immersions $g:T'\to T$ in
$\cS\text{-}\Fan$, in the following sense. Set $d':=\dim\,T'$, and
notice that $d'\leq d$; denote by $\pi_{dd'}:\Q_+^{\oplus
d}\to\Q_+^{\oplus d'}$ the projection on the first $d'$ direct
summands; then the diagram of $T'$-monoids :
\set\begin{equation}\label{eq_compat-coordinates}
{\diagram
(\Q_+^{\oplus d})_{T'(\rho_{T'})} \ar[r]^-{\tilde g^*\theta_T}
\ar[d]_{(\pi_{dd'})_{T'(\rho_{T'})}} & \tilde g^*\cO_{T(\rho_T),\Q}
\ar[d]^{(\log\tilde g)_\Q} \\
(\Q_+^{\oplus d'})_{T'(\rho_{T'})} \ar[r]^-{\theta_{T'}} &
\cO_{T'(\rho_{T'}),\Q}
\enddiagram}
\end{equation}
commutes.

\sset\subsubsection{} Let $f:T'\to T$ be a proper morphism,
with $T'$ locally fine, such that the induced map
$T(\Q_+)\to T'(\Q_+)$ is injective, and let $s\in T$ be any element.
For every $t\in f^{-1}(s)$, we set
$$
G_t:=U(t)(\Q_+)^\gp\cap U(s)(\N)^\gp
\qquad
H_t:=U(t)(\N)^\gp
\qquad
\delta_t:=(G_t:H_t)
$$
and define $\delta(f,s):=\max(\delta_t~|~t\in f^{-1}(s))$.

\begin{lemma}\label{lem_by-dets-I-subdivide}
For every $s\in T$ we have :
\begin{enumerate}
\item
$\delta(f,s)\in\N$.
\item
If $t,t'\in f^{-1}(s)$, and $t$ is a specialization of\/ $t'$ in
$T'$, then $\delta_{t'}\leq\delta_t$.
\item
If\/ $f$ is a rational subdivision, the following conditions
are equivalent :
\begin{enumerate}
\item
$\delta(f,s)=1$.
\item
$U(s)(\N)=\bigcup_{t\in f^{-1}(s)}U(t)(\N)$.
\end{enumerate}
\end{enumerate}
\end{lemma}
\begin{proof}(i): Since $f^{-1}(s)$ is a finite set, the assertion
means that $\delta_t<+\infty$ for every $t\in f^{-1}(s)$. However,
for such a $t$, let $P:=\cO_{T,s}$ and $Q:=\cO_{T',t}$; then
$$
U(t)(\Q_+)^\gp=\Hom_\Z(Q^\gp,\Q)
\qquad\text{and}\qquad
U(t)(\N)^\gp=\Hom_\Z(Q^\gp,\Z)
$$
(remark \ref{rem_reflex-dual}(i) and proposition
\ref{prop_reflex-dual}(iii)). By proposition
\ref{prop_val-crit-fans} we know that the map $(\log
f)^\gp_t\otimes_\Z\Q:P^\gp\otimes_\Z\Q\to Q^\gp\otimes_\Z\Q$ is
surjective. Since $Q$ is finitely generated, it follows that the
image $Q'$ of $P^\gp$ in $Q^\gp$ is a subgroup of finite index, and
then it is easily seen that $\delta_t=(Q^\gp:Q')$.

(ii): Under the stated assumptions we have :
$$
G_{t'}=G_t\cap U(t')(\Q_+)^\gp \qquad H_{t'}=H_t\cap U(t')(\Q_+)^\gp
$$
whence the contention.

(iii): Assume that (b) holds, and let $\phi\in U(t)(\Q_+)^\gp\cap
U(s)(\N)^\gp$ for some $t\in f^{-1}(s)$. Pick any element
$\phi_0\in U(t)(\N)$ which does not lie on any proper face of
$\cO_{T',t}^\vee$; we may then find an integer $a>0$ large enough, so
that $\phi+a\phi_0\in U(t)(\Q_+)$. Set $\phi_1:=\phi+a\phi_0$ and
$\phi_2:=a\phi_0$; then $\phi_1,\phi_2\in U(t)(\Q_+)\cap
U(s)(\N)=U(t)(\N)$, hence $\phi=\phi_1-\phi_2\in U(t)(\N)^\gp$.
Since $\phi$ is arbitrary, we see that $\delta_t=1$, whence (a).

Conversely, suppose that (a) holds, and let $\phi\in U(s)(\N)$; then
there exists $t\in f^{-1}(s)$ such that $\phi\in U(t)(\Q_+)$.
Thus, $\phi\in U(t)(\N)^\gp\cap U(t)(\Q_+)=U(t)(\N)$, whence (b).
\end{proof}

\begin{theorem}\label{th_Saint-Donat}
Every locally fine and saturated fan $T$ admits an integral, proper,
simplicial subdivision $f:T'\to T$, whose restriction
$f^{-1}T_\mathrm{sim}\to T_\mathrm{sim}$ is an isomorphism of fans.
\end{theorem}
\begin{proof} Let $T$ be such a fan. By induction on $h\in\N$, we
shall construct a system of integral, proper simplicial
subdivisions $f_h:S(h)\to T_h$ of the open subsets $T_h$ (notation
of \eqref{subsec_height-in-T}), such that, for every $h\in\N$, the
restriction $f_{h+1}^{-1}(T_h)\to T_h$ of $f_{h+1}$ is isomorphic
to $f_h$, and such that
$f_h^{-1}T_{h,\mathrm{sim}}\to T_{h,\mathrm{sim}}$ is an isomorphism.
Then, the colimit of the morphisms $f_h$ will be the sought
subdivision of $T$.

For $h\leq 1$, we may take $S(h):=T_h$.

Next, suppose that $h>1$, and that $f_{h-1}:S(h-1)\to T_{h-1}$ has
already been given; as a first step, we shall exhibit a rational,
proper simplicial subdivision of $T_h$.
Indeed, for every $t\in T_h\setminus T_{h-1}$, choose a $\N$-point
$\phi_t\in U(t)(\N)$ in the following way. If $t\in T_\mathrm{sim}$,
then let $\phi_t$ be the (unique) generator of an arbitrarily chosen
one-dimensional face of $U(t)(\N)$; and otherwise take any point
$\phi_t$ which does not lie on any proper face of $U(t)(\N)$.

With these choices, notice that
$U(\phi_t)\times(\Spec\,\N)^\sharp=U(t)$ in case
$t\in T_\mathrm{sim}$, and otherwise $U(\phi_t)=U(t)_{h-1}$
(notation of \eqref{subsec_std-subdivision}).
By lemma \ref{lem_star-subdivide}, we obtain corresponding rational
subdivisions $U(\phi_t)\times(\Spec\,\N)^\sharp\to U(t)$ of $U(t)$
centered at $\phi_t$. Notice that if $t$ lies in the simplicial
locus, this subdivision is an isomorphism, and in any case, it
restricts to an isomorphism on the preimage of $U(t)_{h-1}$.

By composing with the restriction of
$f_{h-1}\times(\Spec\,\N)^\sharp$, we get a rational subdivision:
$$
g_t:
T'_t:=f^{-1}_{h-1}U(\phi_t)\times(\Spec\,\N)^\sharp\to U(t)
$$
whose restriction to the preimage of $U(t)_{h-1}$ is an
isomorphism. Moreover, $g_t$ is an isomorphism if $t$ lies
in $T_\mathrm{sim}$. Also notice that $g_t$ is simplicial,
since the same holds for $f_{h-1}$.

If $t,t'$ are any two distinct points of $T$ of height $h$, we
deduce an isomorphism
$$
g_t^{-1}(U(t)\cap U(t'))\isom g^{-1}_{t'}(U(t)\cap U(t'))
$$
hence we may glue the fans $T'_t$ and the morphisms $g_t$
along these isomorphism, to obtain the sought simplicial rational
subdivision $g:T'\to T_h$.

For the next step, we shall refine $g$ locally at every point
$s$ of height $h$; {\em i.e.} for such $s$, we shall find an integral,
proper simplicial subdivision $f_s:T''_s\to U(s)$, whose restriction to
$U(s)_{h-1}$ agrees with $g$, hence with $f_{h-1}$.
Once this is accomplished, we shall be able to build the sought
subdivision $f_h$ by gluing the morphisms $f_s$ and $f_{h-1}$
along the open subsets $U(s)_{h-1}$.

Of course, if $s$ lies in the simplicial locus of $T$, we
will just take for $f_s$ the restriction of $g$, which by
construction is already an isomorphism.

Henceforth, we may assume that $T=U(s)$ is an affine fan of
dimension $h$ with $s\notin T_\mathrm{sim}$, and $g:T'\to T$
is a given proper rational simplicial subdivision, whose
restriction to $g^{-1}T_{h-1}$ is an integral subdivision.
We wish to apply the criterion of lemma
\ref{lem_by-dets-I-subdivide}(iii), which shows that $g$ is an integral
subdivision if and only if $\delta(g,s)=1$. Thus, let $t_1,\dots,t_k$
be the points of $T'$ such that $\delta_{t_i}=\delta(g,s)$ for
every $i=1,\dots,k$. Since $\delta(g,s)$ is anyway a positive
integer (lemma \ref{lem_by-dets-I-subdivide}(i)), a simple
descending induction reduces to the following :

\begin{claim} Given $g$ as above, we may find a proper rational
simplicial subdivision $g':T''\to T$ such that the following holds :
\begin{enumerate}
\item
The restriction of $g'$ to $g^{\prime -1}T_{h-1}$
is isomorphic to the restriction of $g$.
\item
Let $t'_1,\dots,t'_{k'}\in T''$ be the points such that
$\delta_{t'_i}=\delta(g',s)$ for every $i=1,\dots,k'$. We have
$\delta(g',s)\leq\delta(g,s)$, and if $\delta(g',s)=\delta(g,s)$
then $k'<k$.
\end{enumerate}
\end{claim}
\begin{pfclaim}[] Set $t:=t_1$; by definition, there exists
$\phi'\in U(s)(\N)$ which lies in $U(t)(\Q_+)\setminus U(t)(\N)$;
this means that there exists a morphism $\phi$ fitting into
a commutative diagram :
$$
\xymatrix{ \cO_{T,s} \ar[rr]^-{(\log g)_t} \ar[d]_{\phi'} & &
\cO_{T',t} \ar[d]^\phi \\ \N \ar[rr] & & \Q_+.}
$$
Say that $\cO_{T',t}\simeq\N^{\oplus r}$, and let
$(\pi_1,\dots,\pi_r)$ be the (essentially unique) basis of
$\cO_{T'(h),t}^\vee$; then $\phi=a_1\pi_1+\cdots a_r\pi_r$, for some
$a_1,\dots,a_r\geq 0$, and after subtracting some positive integer
multiple of $\pi$, we may assume that $0\leq a_i<1$ for every
$i=1,\dots,r$. Moreover, the coefficients are all strictly positive
if and only if $\phi$ is a local morphism; more generally, we let
$t'$ be the unique generization of $t$ such that $\phi$ factors
through a local morphism $\cO_{T',t'}\to\N$. Set $e:=\hgt(t')$,
and denote by $\pi_1,\dots,\pi_e$ the basis of $\cO_{T',t'}^\vee$, so
that :
\set\begin{equation}\label{eq_unique-basis}
\phi=b_1\pi_1+\cdots+b_e\pi_e
\end{equation}
for unique rational coefficients $b_1,\dots,b_e$ such that $0<b_i<1$
for every $i=1,\dots,e$.

Denote by $Z\subset T'$ the topological closure of $\{t'\}$; for
every $u\in Z\cap T'$, the morphism $\phi$ factors through a
morphism $\phi_u:\cO_{T,u}\to\Q_+$,
and we may therefore consider the subdivision of $U(u)$ centered at
$\phi_u$ as in \eqref{subsec_std-subdivision}, which fits into a
commutative diagram of fans :
$$
\xymatrix{
(U(u)\!\setminus\!Z)\times(\Spec\,\Q_+)^\sharp \ar[r] \ar[d] &
(U(u)\!\setminus\!Z)\times(\Spec\,\N)^\sharp \ar[d]^{\psi_u} \\
U(u) \ar[r] & U'(u):=(\Spec\,\phi_u^{-1}\N)^\sharp.}
$$
We complete the family
$(\psi_u~|~u\in Z)$, by letting $U'(u):=U(u)$ and
$\psi_u:=\one_{U(u)}$ for every $u\in T'\setminus Z$.
Notice then, that the topological spaces underlying $U(u)$ and
$U'(u)$ agree for every $u\in T'$, and for every $u_1,u_2\in T'$,
the restrictions of $\psi_{u_1}$ and $\psi_{u_2}$ :
$$
\psi_{u_i}^{-1}(U(u_1)\cap U(u_2))\to U(u_1)\cap U(u_2)
\qquad (i=1,2)
$$
are isomorphic. Furthermore, by construction, each restriction
$U(u)\to T$ of $g$ factors uniquely through a morphism
$\beta_u:U'(u)\to T$, hence the family
$(\beta_u\circ\psi_u~|~u\in T')$ glues to a well defined
morphism of fans $g':T''\to T$. By a direct inspection, it
is easily seen that $g'$ is a proper rational simplicial
subdivision which fulfills condition (i) of the claim.

Moreover, the map of topological spaces underlying $g'$
factors naturally through a continuous map $p:T''\to T'$,
so that $g'=g\circ p$. The restriction
$V:=p^{-1}(T'\setminus\!Z)\to T$ of $g'$ is isomorphic to
the restriction $T'\setminus\!Z\to T$ of $g$, hence :
$$
\delta_t=\delta_{p(t)} \qquad \text{for every $t\in
p^{-1}(T'\setminus\!Z)$.}
$$
It follows that if $k>1$ and $T'\setminus\!Z$ contains at
least one of the points $t_2,\dots, t_k$, then
$\delta(g',s)\geq\delta(g,s)$. On the other hand,
$p^{-1}(T'\setminus\!Z)$ contains at most $k-1$ points $u$
of $T''$ such that $\delta_u=\delta(g,s)$, and for the remaining
points $u'\in p^{-1}(T'\setminus\!Z)$ we have
$\delta_{u'}<\delta(g,s)$. Since obviously
$$
\delta(g',s)=
\max(\delta(g'_{|p^{-1}(T'\setminus\!Z)},s),\delta(g'_{|p^{-1}Z},s))
$$
we see that condition (ii) holds provided we show :
\set\begin{equation}\label{eq_provided-we-show}
\delta(g'_{|p^{-1}Z},s)=\max(b_1,\dots,b_e)\cdot\delta(g,s).
\end{equation}
Hence, let us fix $u\in Z$, and let
$(v,x)\in(U(u)\!\setminus\!Z)\times(\Spec\,\N)^\sharp$ be any point
(see \eqref{prod_of-fine-fans}); if $x=\emptyset$, then
$(v,x)\notin p^{-1}Z$, so it suffices to consider the points
of the form $(v,\fm)$ (where $\fm\subset\N$ is the maximal ideal).
Moreover, say that $\hgt(u)=d$; in view of lemma
\ref{lem_by-dets-I-subdivide}(ii), it suffices to consider the
points $(v,\fm)$ such that $\hgt(v)=d-1$. There are
exactly $e$ such points, namely the prime ideals
$v_i:=(\pi_i\circ\sigma)^{-1}\fm$, where $\pi_1,\dots,\pi_e$
are as in \eqref{eq_unique-basis},  and
$\sigma:\cO_{T',u}\to\cO_{T',t'}$ is the specialization map.

In order to estimate $\delta_{(v,\fm)}$ for some $v:=v_i$, we
look at the transpose of the map
$$
(\log g')^\gp_{(v,\fm)}:\cO_{T,s}^\gp\to\cO_{T',v}^\gp\times\Z
\quad : \quad
z\mapsto((\log g)_v^\gp(z),\phi^\gp(z)).
$$
Let $(p_1,\dots,p_d)$ be the basis of $\cO_{T,u}^\vee$, ordered
in such a way that $p_i=\pi_i\circ\sigma$ for every $i=1,\dots,e$.
By a little abuse of notation, we may then denote $(p_j~|~j\neq i)$
the basis of $\cO_{T',v}^\vee$, so that the dual group
$(\cO_{T',v}^\gp\times\Z)^\vee$ admits the basis
$\{p_j^\gp~|~j\neq i\}\cup\{q\}$, where $q$ is the natural
projection onto $\Z$.
Set as well $p'_i:=p_i\circ(\log g)_u^\gp$ for every $i=1,\dots,d$.
With this notation, the above transpose is the group homomorphism
given by the rule :
$$
p_j\mapsto p'_i \quad \text{for $j\neq i$} \qquad \text{and} \qquad
q\mapsto b_1p'_1+\cdots+b_ep'_e
$$
from which we deduce easily that $\delta_{(v,\fm)}=b_i\cdot\delta_u$,
whence \eqref{eq_provided-we-show}.
\end{pfclaim}
\end{proof}

\begin{proposition}\label{prop_flattening}
Let $(\Gamma,+,0)$ be a fine monoid, $M$ a fine $\Gamma$-graded
monoid. Then :
\begin{enumerate}
\item
There exists a finite set of generators $C:=\{\gamma_1,\dots,\gamma_k\}$
of\/ $\Gamma$, with the following property. For every $\gamma\in\Gamma$,
we may find $a_1,\dots,a_k\in\N$ such that:
$$
\gamma=a_1\gamma_1+\cdots+a_k\gamma_k
\quad\text{and}\quad
M_\gamma=M^{a_1}_{\gamma_1}\cdots M^{a_k}_{\gamma_k}.
$$
\item
There exists a subgroup $H\subset\Gamma^\gp$ of finite index,
such that :
$$
M_{a\gamma}=M_\gamma^a
\qquad
\text{for every $\gamma\in H\cap\Gamma$ and every integer $a>0$}.
$$
\end{enumerate}
(In {\em (i)} and {\em (ii)} we use the multiplication law of\/
$\cP(M)$, as in \eqref{subsec_toric}).
\end{proposition}
\begin{proof} Obviously we may assume that $M$ maps surjectively
onto $\Gamma$, in which case $G:=\Gamma^\gp$ is finitely generated,
and its image $G'$ into $G_\R:=G\otimes_\Z\R$ is a free abelian
group of finite rank. The same holds as well for the image $L$
of $\log M^\gp$ in $M^\gp_\R:=\log M^\gp\otimes_\Z\R$.
Let $p:G\to G'$ be the natural projection. According to
proposition \ref{prop_Gordon}(iii), we have :
\set\begin{equation}\label{eq_indeed}
M_\Q=M_\R\cap M^\gp_\Q
\end{equation}
(notation of \eqref{subsec_from-con-to-mon}). Let
$f^\gp_\R:M^\gp_\R\to G_\R$ be the induced $\R$-linear map, and
denote by $f_\R:M_\R\to G_\R$ the restriction of $f^\gp_\R$. By
proposition \ref{prop_divide-et-imp}(ii), we may find a $G'$-rational
subdivision $\Delta$ of $f_\R(M_\R)$ such that :
\set\begin{equation}\label{eq_good-fan}
f^{-1}_\R(a+b)=f^{-1}_\R(a)+f^{-1}_\R(b)
\end{equation}
for every $\sigma\in\Delta$ and every $a,b\in\sigma$. After
choosing a refinement, we may assume that $\Delta$ is a simplicial
fan (theorem \ref{th_Saint-Donat}). Let $\tau\in\Delta$ be any
cone; by proposition \ref{prop_Gordon}(i), the monoid
$N:=\tau\cap G'$ is finitely generated, and then the same holds
for $M\times_{G'}N$, by corollary \ref{cor_fibres-are-fg}. However,
the latter is just $M':=\bigoplus_{\gamma\in p^{-1}N}M_\gamma$
(lemma \ref{lem_forget-me-not}(iii)). Set $\Gamma':=\Gamma\cap p^{-1}N$.

\begin{claim}\label{cl_reduce-to-simpl}
$p(\Gamma')$ generates $\tau$.
\end{claim}
\begin{pfclaim}
Clearly the $\Gamma$-grading of $M$ induces a surjection
$M_\Q\to\Gamma_{\!\Q}$ (notation of \eqref{subsec_from-con-to-mon}).
By proposition \ref{prop_Gordon}(iii), we have
$\Gamma_{\!\Q}=f_\R(M_\R)\cap G_\Q$, hence
$N\subset\tau\cap G_\Q\subset\Gamma_{\!\Q}$. Thus, for every $n\in N$
we may find $a\in\N$ such that $a\cdot n\in p(\Gamma)$, hence
$a\cdot n\in p(\Gamma')$. On the other hand, $N$ generates $\tau$,
since the latter is $G'$-rational. The claim follows.
\end{pfclaim}

In view of claim \ref{cl_reduce-to-simpl}, we may replace $M$ by
$M'$, and $\Gamma$ by $\Gamma'$, which allows to assume that
$f_\R(M_\R)$ is a simplicial cone, and \eqref{eq_good-fan} holds
for every $a,b\in f_\R(M_\R)$. Next, let
$S:=\{e_1,\dots,e_n\}\subset p(\Gamma)$ be a set of generators of
the cone $f_\R(M_\R)$; from the discussion in
\eqref{subsec_extremal} we see that, up to replacing $S$ by a
subset, the rays $\R_+\cdot e_i$ (with $i=1,\dots,n$) are precisely
the extremal rays of $f_\R(M_\R)$, especially, the vectors
$e_1,\dots,e_n$ are $\R$-linearly independent. Choose
$g_1,\dots,g_n\in\Gamma$ such that $p(g_i)=e_i$ for every
$i=1,\dots,n$. According to corollary \ref{cor_no-fibres-here},
there exist finite subsets $\Sigma_1\dots,\Sigma_n\subset M$,
such that :
$$
M_{g_i}=M_0\cdot\Sigma_i
\qquad
\text{for every $i=1,\dots,n$}.
$$

\begin{claim}\label{cl_cucu}
$(M_0)_\Q=f^{-1}_\R(0)\cap M^\gp_\Q$.
\end{claim}
\begin{pfclaim} To begin with, we may write
$f^{-1}_\R(0)=(f^\gp_\R)^{-1}(0)\cap M_\R$, hence $f^{-1}_\R(0)\cap
M^\gp_\Q=(f^\gp_\Q)^{-1}(0)\cap M_\Q$, by \eqref{eq_indeed}. Now,
suppose $x\in (f^\gp_\Q)^{-1}(0)\cap M_\Q$; then we may find an
integer $a>0$ such that $ax=m\otimes 1$ for some $m\in\log M$.
Say that $m\in\log M_\gamma$; then $f^\gp_\Q(\gamma)=0$, therefore
$\gamma$ is a torsion element of $G'$, and consequently $bm\in\log M_0$
for an integer $b>0$ large enough. We conclude that
$x=(bm)\otimes(ba)^{-1}\in(M_0)_\Q$, as claimed.
\end{pfclaim}

On the other hand, since $\R_+e_i$ is a $G'$-rational polyhedral
cone, $f^{-1}_\R(\R_+e_i)$ is an $L$-rational polyhedral cone
(proposition \ref{prop_was-part-of-Gordon}(ii,iii)), hence it
admits a finite set of generators $S_i\subset L$. Up to replacing
the elements of $S_i$ by some positive rational multiples, we may
assume that $f_\R(s)$ is either $0$ or $e_i$, for every $s\in S_i$.
In this case, it follows easily that
\set\begin{equation}\label{eq_conv-hull}
f^{-1}_\R(e_i)=f^{-1}_\R(0)+T_i
\end{equation}
where :
$$
T_i:= \biggl\{\sum_{s\in S_i}t_s\cdot s~|~t_s\in\R_+,\sum_{s\in
S_i}t_s=1\biggr\}
$$
is the convex hull of $S_i$ (and as usual, the addition of sets in
\eqref{eq_conv-hull} refers to the addition law of $\cP(M^\gp_\R)$,
see \eqref{subsec_toric}). Next, notice that $S_i\subset M_\Q$, by
\eqref{eq_indeed}; thus, we may find an integer $a>0$ such that
$a\cdot s$ lies in the image of $\log M$, for every $s\in S_i$.
After replacing $e_i$ by $a\cdot e_i$ and $S_i$ by $\{a\cdot
s~|~s\in S_i\}$, we may then achieve that \eqref{eq_conv-hull}
holds, and furthermore $S_i$ lies in the image of $\log M$,
therefore in the image of $\log M_{g_i}$. It follows easily that
\eqref{eq_conv-hull} still holds with $S_i$ replaced by the set
$\Sigma_i\otimes 1:=\{m\otimes 1~|~m\in\Sigma_i\}$. Let
$\Sigma_0\subset\log M$ be a finite set of generators for the monoid
$M_0$; claim \ref{cl_cucu} implies that $\Sigma_0$ is also a
set of generators for the $L$-rational polyhedral cone
$f_\R^{-1}(0)$. Let $P\subset M$ be the submonoid generated by
$\Sigma:=\Sigma_0\cup\Sigma_1\cup\cdots\cup\Sigma_n$, and
$\Delta\subset\Gamma$ the submonoid generated by $g_1,\dots,g_n$;
clearly the $\Gamma$-grading of $M$ restricts to a $\Delta$-grading
on $P$. Notice that $g_1\otimes 1,\dots,g_n\otimes 1$ are linearly
independent in $G_\Q$, since the same holds for $e_1,\dots,e_n$;
especially, $\Delta\simeq\N^{\oplus n}$.

\begin{claim}\label{cl_let-A}
(i)\ \ The set $\Sigma\otimes 1:=\{s\otimes 1~|~s\in\Sigma\}$
generates the cone $M_\R$.
\begin{enumerate}
\addenu
\item
$P_{a+b}=P_a\cdot P_b$ for every $a,b\in\Delta$.
\item
There exists a finite set $A\subset M$ such that $M=A\cdot P$.
\end{enumerate}
\end{claim}
\begin{pfclaim} By \eqref{eq_conv-hull} and the foregoing discussion,
we know that $(\Sigma_0\cup\Sigma_i)\otimes 1$ generate
$f^{-1}_\R(\R_+e_i)$, for every $i=1,\dots,n$. Since the
additivity property \eqref{eq_good-fan} holds for every $a,b\in
f_\R(M_\R)$, assertion (i) follows. (ii) is a straightforward
consequence of the definitions. Next, from (i) and proposition
\ref{prop_Gordon}(iii), we deduce that $M_\Q=P_\Q$. Thus, let
$m_1,\dots,m_r$ be a system of generators for the monoid $M$;
it follows that there are integers $k_1,\dots,k_r>0$, such that
$m_1^{k_1},\dots,m_r^{k_r}\in P$, and therefore the subset
$A:=
\{\prod_{i=1}^rm_i^{t_i}~|~\text{$0\leq t_i<k_i$ for every $i\leq r$}\}$
fulfills the condition of (iii).
\end{pfclaim}

We introduce a partial ordering on $G$, by declaring that $a\leq b$
for two elements $a,b\in G$, if and only if $b-a\in\Delta$. Now, let
$a\in G$ be any element; we set
$$
G(a):=\{g\in G~|~g\leq a\}.
$$
The subset $G(a)$ inherits a partial ordering from $G$, and $a$ is
the maximum of the elements of $G(a)$; moreover, notice that every
finite subset $S\subset G(a)$ admits a supremum $\sup S\in G(a)$.
Indeed, it suffices to show the assertion for a set of two elements
$S=\{b_1,b_2\}$; we may then write $a-b_i=\sum_{j=1}^nk_{ij}g_j$ for
certain $k_{ij}\in\N$, and then
$\sup(b_1,b_2)=a-\sum_{j=1}^n\min(k_{1j},k_{2j})\cdot g_j$.

Let $A$ be as in claim \ref{cl_let-A}(iii), and denote by
$B\subset\Gamma$ the image of $A$; then
$$
M=M_0\cdot M_B
\qquad
\text{where\ \ $M_B:=\bigoplus_{b\in B}\log M_b$.}
$$
For every $a\in G$, let also $B(a):=B\cap G(a)$; invoking several
times claim \ref{cl_let-A}(ii), we get :
\begin{align*}
M_a & =\bigcup_{b\in B(a)}M_b\cdot P_{a-b} \\
          & =\bigcup_{b\in B(a)}M_b\cdot P_{\sup B(a)-b}\cdot
             P_{a-\sup B(a)} \\
          & \subset M_{\sup B(a)}\cdot P_{a-\sup B(a)}.
\end{align*}
Finally, say that $a-\sup B(a)=\sum^n_{i=1}t_ig_i$ for certain
$t_1,\dots,t_n\in\N$; applying once more claim \ref{cl_let-A}(ii),
we conclude that :
$$
M_a\subset M_{\sup B(a)}\cdot\prod^n_{i=1}M_{g_i}^{t_i}.
$$
The converse inclusion is clear, and therefore the set
$C:=\{g_1,\dots,g_n\}\cup\{\sup B(a)~|~a\in G\}$ fulfills condition
(i) of the proposition.

(ii): For $h\in\Delta^\gp$, say $h=\sum_{i=1}^n a_ig_i$, with integers
$a_1,\dots,a_n$, we let $|h|:=\sum_{i=1}^n|a_i|g_i\in\Delta$. Choose
any positive integer $\alpha$ such that :
\set\begin{equation}\label{eq_choose-alpha}
|b|\leq\alpha\cdot\sum_{i=1}^ng_i
\qquad
\text{for every $b\in B\cap\Delta^\gp$}
\end{equation}
and let $H\subset\Delta^\gp$ be the subgroup generated by
$\alpha g_1,\dots,\alpha g_n$.

\begin{claim}\label{cl_all-the-same}
$B(h)=B(kh)$ for every $h\in H$ and every integer $k>0$.
\end{claim}
\begin{pfclaim} Let $h:=\sum_{i=1}^n\alpha_i g_i\in H$, and
suppose that $b\in B(kh)$ for some $k>0$; therefore $kh-b\in\Delta$,
hence $b\in\Delta^\gp$, and we can write $b=\sum_{i=1}^n\beta_i g_i$
for integers $\beta_1,\dots,\beta_n$, such that
$k\alpha_i-\beta_i\geq 0$ for every $i=1,\dots,n$. In this case,
\eqref{eq_choose-alpha} implies that $\alpha_i\geq 0$ for every
$i\leq n$, and $\alpha_i\geq\alpha\geq\beta_i$ whenever $\beta_i>0$.
It follows easily that $k'h-b\in\Delta$ for every integer $k'>0$,
whence the claim.
\end{pfclaim}

Using claims \ref{cl_let-A}(ii) and \ref{cl_all-the-same}, and
arguing as in the foregoing, we may compute :
\begin{align*}
M_{ah} & =M_{\sup B(ah)}\cdot P_{ah-\sup B(ah)} \\
& =M_{\sup B(h)}\cdot P_{ah-\sup B(h)} \\
& =M_{\sup B(h)}\cdot P_{h-\sup B(h)}\cdot P_h^{a-1} \\
& \subset M_h^a.
\end{align*}
The converse inclusion is clear, so (ii) holds.
\end{proof}

\sset\subsubsection{} Let $M$ be an integral monoid, and $w\in\log
M^\gp$ any element. For $\eps\in\{1,-1\}$ we have a natural
inclusion
$$
j_\eps:\log M\to M(\eps):=\log M+\eps\N w
$$
({\em i.e.} $M(\eps)$ is the submonoid of $M^\gp$ generated by $M$
and $w^\eps$). Let us write $w:=b^{-1}a$ for some $a,b\in M$; then
the induced morphisms of affine schemes
$\iota_\eps:=\Spec\,\Z[j_\eps]$ have a natural geometric
interpretation. Namely, let $f:X\to\Spec\,\Z[M]$ be the blow up of
the ideal $I\subset\Z[M]$ generated by $a$ and $b$; we have
$X=U_1\cup U_{-1}$, where $U_\eps$, for $\eps=\pm 1$, is the largest
open subscheme of $X$ such that $w^\eps\in\cO_{\!X}(U_\eps)$. Then
$\iota_\eps$ is naturally identified with the restriction
$U_\eps\to\Spec\,\Z[M]$ of the blow up $f$. More generally, by
adding to $M$ any finite number of elements of $M^\gp$, we may
construct in a combinatorial fashion, the standard affine charts of
a blow up of an ideal of $\Z[M]$ generated by finitely many elements
of $M$. These considerations explain the significance of the
following {\em flattening theorem} :

\begin{theorem}\label{th_flattening-th}
Let $j:M\to N$ be an inclusion of fine monoids. Then there
exist a finite set\/ $\Sigma\subset\log M^\gp$, and an
integer $k>0$ such that the following holds :
\begin{enumerate}
\item
For every mapping $\eps:\Sigma\to\{\pm 1\}$, the induced inclusion :
$$
M(\eps):=\log M+\sum_{\sigma\in\Sigma}\eps(\sigma)\N\sigma\to
N(\eps):=\log N+\sum_{\sigma\in\Sigma}\eps(\sigma)\N\sigma
$$
is a flat morphism of fine monoids.
\item
Suppose that $j$ is a flat morphism, and let $\iota:M\to M^\sat$ be
the natural inclusion. Denote by $Q$ the push-out of the diagram
$N\xleftarrow{j}M\xrightarrow{\iota\circ\bek_M}M^\sat$ (where
$\bek_M$ is the $k$-Frobenius). Then the natural map
$M^\sat\to Q^\sat$ is flat and saturated.
\end{enumerate}
\end{theorem}
\begin{proof} (i): (Notice that $\log M$, $\log N$ and
$P(\eps):=\sum_{\sigma\in\Sigma}\eps(\sigma)\N\sigma$ may be
regarded as submonoids of $\log N^\gp$, and then the above sum is
taken in the monoid $(\cP(\log N^\gp),+)$ defined as in
\eqref{subsec_toric}.) Set $G:=N^\gp/M^\gp$, and let
$N^\gp=\bigoplus_{\gamma\in G}N^\gp_\gamma$ be the $j$-grading
of $N^\gp$ (remark \ref{rem_why-not}(iii)); notice that this
grading restricts to $j$-gradings for $N$ and $N(\eps)$. Let
also $\Gamma:=\{\gamma\in G~|~N_\gamma\neq\emptyset\}$, and
choose a finite generating set $\{\gamma_1,\dots,\gamma_r\}$
of $\Gamma$ with the properties of proposition
\ref{prop_flattening}(i). According to corollary
\ref{cor_no-fibres-here}, for every $i\leq r$ there exists a finite
subset $\Sigma_i:=\{t_{i1},\dots,t_{in_i}\}\subset M$ such that
$N_{\gamma_i}=N_0\cdot\Sigma_i$. Moreover, $N_0$ is a
finitely generated monoid, by corollary \ref{cor_fibres-are-fg}.
Let $\Sigma_0$ be a finite system of generators for $N_0$, and
for every $i\leq r$ define
$\Sigma_i':=\{t_{ij}-t_{il}~|~1\leq j<l\leq n_i\}$. We claim that
the subset $\Sigma:=\Sigma_0\cup\Sigma'_1\cup\cdots\cup\Sigma'_r$
will do. Indeed, first of all notice that $\Sigma\subset\log M^\gp$.
We shall apply the flatness criterion of remark \ref{rem_why-not}(iv).
Thus, we have to show that, for every $g\in\Gamma$, the
$M(\eps)$-module $N(\eps)_g$ is a filtered union of cosets
$\{m\}+M(\eps)$ (for certain $m\in N_g$).
Hence, let $x_1,x_2\in N(\eps)_g$; we may write $x_i=m_i+p_i$
for some $m_i\in\log N$ and $p_i\in P(\eps)$ ($i=1,2$); since
$\{x_i\}+M(\eps)\subset\{m_i\}+M(\eps)$, we may then assume
that $p_1=p_2=0$, hence $x_1,x_2\in N$. Thus, it suffices to
show that $N_g$ is contained in a filtered union of cosets
as above. However, by assumption there exist $a_1,\dots,a_r\in\N$
such that $g=\sum_{i=1}^r a_i\gamma_i$ and
$N_g=N_{\gamma_1}^{a_1}\cdots N^{a_r}_{\gamma_r}$; we are then
easily reduced to the case where $g=\gamma_i$ for some $i\leq r$.
Therefore, $x_j=\sigma_j+b_j$, where $\sigma_j\in\Sigma_i$ and
$b_j\in N_0$, for $j=1,2$. Notice that $P(\eps)$ contains either
$\sigma_1-\sigma_2$, or $\sigma_2-\sigma_1$ (or both); in the
first occurrence, set $\sigma:=\sigma_2$, and otherwise, let
$\sigma:=\sigma_1$. Likewise, say that $\Sigma_0=\{y_1,\dots,y_n\}$,
so that $b_j=\sum_{s=1}^n a_{js}y_s$ for certain $a_{js}\in\N$
($j=1,2$); for every $s\leq n$, we set $a^*_s:=\min(a_{1s},a_{2s})$
if $y_s\in P(\eps)$, and otherwise we set $a^*_s:=\max(a_{1s},a_{2s})$.
One sees easily that
$x_1,x_2\in\{\sigma+\sum_{s=1}^n a^*_s y_s\}+M(\eps)$, whence the
contention.

(ii): By proposition \ref{prop_flattening}(ii), there exists a
subgroup $H\subset G$ of finite index, such that :
$$
\pi^{-1}(h^n)=\pi^{-1}(h)^n
$$
for every integer $n>0$ and every $h\in H$. Let $k:=(G:H)$, and
define $N'$ as the fibre product in the cartesian diagram :
\set\begin{equation}\label{eq_for-Nprime}
{\diagram N' \ar[r]^-{\pi'} \ar[d]_\mu & G \ar[d]^{\bek_G} \\
           N \ar[r]^-\pi & G
\enddiagram}
\end{equation}
The trivial morphism $\boldsymbol{0}_M:M\to G$ ({\em i.e.} the
unique one that factors through $\{1\}$) and the inclusion $j$
satisfy the identity : $\bek_G\circ\boldsymbol{0}_M=\pi\circ j$,
hence they determine a well-defined map $\phi:M\to N'$.

\begin{claim}\label{cl_phi-fl-sat} $\phi$ is flat and saturated.
\end{claim}
\begin{pfclaim} For the flatness, we shall apply the criterion
of remark \ref{rem_why-not}(iv). First, $\phi$ is injective, since
$\mu\circ\phi=j$. Next, notice that the sequence of abelian groups :
$$
0\to M^\gp\xrightarrow{\phi^\gp} (N')^\gp\xrightarrow{(\pi')^\gp}
G\to 0
$$
is short exact; indeed, this is none else than the pullback
$\cE*\bek_G^\gp$ along the morphism $\bek_G$, of the short exact
sequence
\set\begin{equation}\label{eq_for-E} \cE:=(0\to
M^\gp\xrightarrow{j^\gp} N^\gp\xrightarrow{\pi^\gp} G\to 0).
\end{equation}
It follows that $\phi$ is flat if and only if, for every $x\in
\pi'(N')$, the preimage $(\pi')^{-1}(x)$ is a filtered union of
cosets of the form $\{n\}\cdot\phi(M')$. However, the induced map
$(\pi')^{-1}(g)\to\pi^{-1}(g^k)$ is a bijection for every $g\in G$,
and the flatness of $j$ implies that $\pi^{-1}(x^k)$ is a filtered
union of cosets, whence the contention. Notice also that
$\Img\,\bek_G\subset H$; hence, by the same token, we derive that
$(\pi')^{-1}(g^n)=(\pi')^{-1}(g)^n$ for every $g\in G$, therefore
$\phi$ is quasi-saturated, by proposition \ref{prop_crit-saturated}.
Since we know already that $j$ is integral (theorem
\ref{th_flat-crit-for-mnds}), the claim follows.
\end{pfclaim}

Next, we wish to consider the commutative diagram of monoids :
\set\begin{equation}\label{eq_for-psi}
{\diagram M \ar[r]^-j \ar[d]_{\bek_M} & N \ar[d]^\psi \\
          M \ar[r]^-\phi & N'
\enddiagram}
\end{equation}
such that $\psi$ is the map determined by the pair of morphisms
$(f,\bek_N)$. Let $P$ be the push-out of the maps $j$ and $\bek_M$;
the maps $\phi$ and $\psi$ determine a morphism $\tau:P\to N'$.

\begin{claim}\label{cl_diags-gp}
(i)\ \  The diagram $\eqref{eq_for-Nprime}^\gp$ of associated
abelian groups, is cartesian.
\begin{enumerate}
\addenu
\item
The diagram of abelian groups $\eqref{eq_for-psi}^\gp$ is
cocartesian ({\em i.e.} $\tau^\gp$ is an isomorphism).
\item
There exists a morphism $\lambda:N'\to P$ such that
$\lambda\circ\tau=\bek_P$ and $\tau\circ\lambda=\bek_{N'}$.
\end{enumerate}
\end{claim}
\begin{pfclaim}(i): Suppose that $x\in(N')^\gp$ is any element
such that $(\pi')^\gp(x)=1$ and $\mu^\gp(x)=1$; we may write
$x=b^{-1}a$ for some $a,b\in N'$, and it follows that
$\pi'(a)=\pi'(b)$ and $\mu(a)=\mu(b)$, hence $a=b$ in $N'$, since
the forgetful functor $\Mnd\to\Set$ commutes with fibre products.
Thus $x=1$. On the other hand, suppose that $\pi^\gp(b^{-1}a)=x^k$
for some $a,b\in N$ and $x\in G$; therefore, $\pi(b^{k-1}a)=(bx)^k$,
so there exists $c\in N'$ such that $\mu(c)=b^{k-1}a$ and
$\pi'(c)=bx$. Likewise, there exists $d\in N'$ with $\mu(d)=b^k$ and
$\pi'(d)=b$. Consequently, $(\pi')^\gp(d^{-1}c)=x$ and
$\mu(d^{-1}c)=b^{-1}a$. The assertion follows.

(ii): Quite generally, let $E:=(0\to A\to B\to C\to 0)$ be a short
exact sequence of objects in an abelian category $\cC$, and for
every object $X$ of $\cC$, let $\bek_X:=k\cdot\one_X:X\to X$. Then
one has a natural map of complexes $\alpha_E:E*\bek_C\to E$ (resp.
$\beta_E:E\to\bek_A*E$) from $E$ to the pull-back of $E$ along
$\bek_C$ (resp. from the push-out of $E$ along $\bek_A$ to $E$), and
a natural isomorphism $\omega_E:\bek_A*E\isom E*\bek_C$ in the
category $\mathbf{Ext}_\cC(C,A)$ of extensions of $C$ by $A$ (this
is the category whose objects are all short exact sequences in $\cC$
of the form $0\to A\to X\to C\to 0$, and whose morphisms are the
maps of complexes which are the identity on $A$ and $C$). These maps
are related by the identities :
\set\begin{equation}\label{eq_for-Exts}
\beta_E\circ\alpha_E\circ\omega_E=k\cdot\one_{\bek_A*E}\qquad
\omega_E\circ\beta_E\circ\alpha_E=k\cdot\one_{E*\bek_C}.
\end{equation}
We leave to the reader the construction of $\omega_E$. In the case
at hand, we obtain a natural isomorphism
$\omega_\cE:\bek^\gp_M*\cE\isom\cE*\bek_G$, where $\cE$ is the short
exact sequence of \eqref{eq_for-E}. An inspection of the
construction shows that the map $\tau^\gp$ is precisely the
isomorphism defined by $\omega_\cE$.

(iii): Let $\mu':N\to P$ be the natural map, and set
$\lambda:=\mu'\circ\mu$. By inspecting the constructions, one checks
easily that $\lambda^\gp$ is the map defined by
$\beta_\cE\circ\alpha_\cE$. Then the assertion follows from
\eqref{eq_for-Exts}.
\end{pfclaim}

Let $P'$ be the push-out of the diagram
$N'\xleftarrow{\phi}M\xrightarrow{\iota}M^\sat$; from claim
\ref{cl_phi-fl-sat}, lemma \ref{lem_obvious-int}(i) and corollary
\ref{cor_exact}(iii), we deduce that the natural map $M^\sat\to P'$
is flat and saturated, hence $P'$ is saturated. On the other hand,
directly from the definitions we get a cocartesian diagram
\set\begin{equation}\label{eq_directly-ondef}
{\diagram P \ar[r] \ar[d]_\tau & Q \ar[d]^{\tau'} \\
           N' \ar[r] & P'.
\enddiagram}\end{equation}
The induced diagram $\eqref{eq_directly-ondef}^\sat$ of saturated
monoids is still cocartesian (remark \ref{rem_satura}(v)); however,
claim \ref{cl_diags-gp}(iii) implies easily that $\tau^\sat$ is an
isomorphism, therefore the same holds for $(\tau')^\sat$, and
assertion (ii) follows.
\end{proof}

\section{Homological algebra}
\label{chap_hom-algebra}
This chapter lays the foundations of homological and
homotopical algebra that shall be needed in the rest
of the treatise. For some standard facts that we do
not repeat here, we use the book \cite{We} as a quick
reference.

\subsection{Complexes in an additive category}
\label{sec_brutal-truncate}
Let $\cA$ be any additive category. We denote by $\sC(\cA)$
the category of {\em (cochain) complexes} of objects of $\cA$.
Hence, an object of $\sC(\cA)$ is a pair
$$
K^\bullet:=(K^\bullet,d^\bullet_K)
$$
consisting of a system of objects $(K^n~|~n\in\Z)$ and
morphisms $(d^n_K:K^n\to K^{n+1}~|~n\in\Z)$ of $\cA$, called
the {\em differentials} of the complex $K^\bullet$, such that
$$
d^{n+1}_K\circ d^n_K=0
\qquad
\text{for every $n\in\Z$}.
$$
The morphisms $\phi^\bullet:K^\bullet\to L^\bullet$ in $\sC(\cA)$
are the systems of morphisms $(\phi^n:K^n\to L^n~|~n\in\Z)$ of
$\cA$ such that
$$
\phi^{n+1}\circ d_K^n=d^{n+1}_L\circ\phi^n
\qquad
\text{for every $n\in\Z$}.
$$
We will usually omit the subscript when referring to the
differentials of a complex, unless there is a danger of
confusion. 
Also, let $I\subset\Z$ be any (bounded or unbounded) interval,
{\em i.e.} $I$ is either of the form $\Z\cap[a,+\infty[$, or
$\Z\cap]-\infty,b]$ (for some $a,b\in\Z$), or the intersection
of any of these two. We shall denote by
$$
\sC^I(\cC)
$$
the full subcategory of $\sC(\cC)$ whose objects are the
complexes $K^\bullet$ such that $K^i=0$ whenever $i\notin I$
(where we fix a zero object $0$ of $\cA$).
For instance, if $I=\Z\cap[a,+\infty[$, then $\sC^I(\cC)$
is also denoted $\sC^{\geq a}(\cC)$, and likewise for the
case of an upper bounded interval. Moreover, we set
$$
\sC^-(\cA):=\bigcup_{n\in\Z}\sC^{\leq n}(\cA)
\qquad
\sC^+(\cA):=\bigcup_{n\in\Z}\sC^{\geq n}(\cA)
\qquad
\sC^b(\cA):=\bigcup_{n\in\Z}\sC^{[-n,n]}(\cA)
$$
so $\sC^-(\cA)$ (resp. $\sC^+(\cA)$, resp. $\sC^b(\cA)$)
is the full subcategory of $\sC(\cA)$ whose objects are
the {\em bounded above} (resp. {\em bounded below},
resp. {\em bounded}) complexes of $\cA$.

\sset\subsubsection{}\label{subsec_brutal-truncate}
For every $a\in\Z$, the inclusion functor
\set\begin{equation}\label{eq_bounded-cpx}
\sC^{\geq a}(\cA)\to\sC(\cA)
\qquad
\text{(resp.\ $\sC^{\leq a}(\cA)\to\sC(\cA)$\ )}
\end{equation}
admits a right (resp. left) adjoint called the
{\em brutal truncation functor}, and denoted
$$
t^{\geq a}:\sC(\cA)\to\sC^{\geq a}(\cA)
\qquad
\text{(resp.\ $t^{\leq a}:\sC(\cA)\to\sC^{\leq a}(\cA)$\ )}.
$$
Namely, for any complex $K^\bullet$, we let $t^{\geq a}(K^\bullet)$
be the unique object of $\sC^{\geq a}(\cA)$ that agrees with
$K^\bullet$ in all degrees $\geq a$, and with the same
differentials as $K^\bullet$, in this range of degrees (and
likewise for $t^{\leq a}(K^\bullet)$).
If $\phi^\bullet:K^\bullet\to L^\bullet$ is any morphism in
$\sC(\cA)$, then $t^{\geq a}\phi^\bullet$ is the unique
morphism $t^{\geq a}K^\bullet\to t^{\geq a}L^\bullet$ which
agrees with $\phi^\bullet$ in all degrees $\geq a$ (and
likewise for $t^{\leq a}\phi^\bullet$).

Moreover, if $\cA$ is an abelian category, \eqref{eq_bounded-cpx}
also admits a left (resp. a right) adjoint
$$
\tau^{\geq a}:\sC(\cA)\to\sC^{\geq a}(\cA)
\qquad
\text{(resp.\ $\tau^{\leq a}:\sC(\cA)\to\sC^{\leq a}(\cA)$\ )}
$$
called the {\em normalized truncation functor} (or just the
{\em truncation functor}). Namely, we fix representatives
for all kernels and cokernels of $\cA$, and for any complex
$K^\bullet$, we let $\tau^{\geq a}(K^\bullet)$ be the unique
object of $\sC^{\geq a}(\cA)$ whose term in degree $a$
equals $\Coker\,d^{a-1}_K$, and which agrees with $K^\bullet$
in all degrees $>a$, with the same differentials as
$K^\bullet$, in this range of degrees; the differential
in degree $a$ is the unique morphism $\Coker\,d^{a-1}_K\to K^{a+1}$
whose composition with the natural map $K^a\to\Coker\,d^{a-1}_K$
equals $d^{a-1}_K$. If $\phi^\bullet:K^\bullet\to L^\bullet$ is
any morphism in $\sC(\cA)$, then $\tau^{\geq a}\phi^\bullet$ is
the unique morphism $\tau^{\geq a}K^\bullet\to\tau^{\geq a}L^\bullet$
which agrees with $\phi^\bullet$ in all degrees $>a$,
and which, in degree $a$, is the morphism
$\Coker\,d_K^a\to\Coker\,d_L^a$ induced by $\phi^a$.

Likewise, $\tau^{\leq a}(K^\bullet)$ is the unique complex
of $\sC^{\leq a}(\cA)$ whose term in degree $a$ equals
$\Ker\,d^a_K$, and which agrees with $K^\bullet$
in all degrees $<a$, with the same differentials as
$K^\bullet$, in all degrees $<a-1$; the differential
in degree $a-1$ is the restriction of $d^{a-1}_K$. Again,
for $\phi^\bullet$ as above, $\tau^{\leq a}\phi^\bullet$
agrees with $\phi^\bullet$ in degrees $<a$, and
is the restriction $\Ker\,d^a_K\to\Ker\,d_L^a$
of $\phi^a$ in degree $a$.

\sset\subsubsection{}\label{subsec_shift}
There is an obvious functor
$$
\cA\to\sC(\cA)
\quad :\quad
A\mapsto A[0]
$$
that sends any object $A$ of $\cA$ to the complex with $A$
{\em placed in degree zero}, {\em i.e.} such that $A[0]^i$
equals $A$ if $i=0$, and equals $0$ otherwise (clearly, there
is a unique such complex). On the other hand, the
{\em shift operator} is the functor
$$
\sC(\cA)\to\sC(\cA)
\quad :\quad K^\bullet\to K^\bullet[1]
$$
given by the rule :
$$
K^\bullet[1]^n:=K^{n+1}
\qquad
d_{K[1]}^n:=-d^{n+1}_K
\qquad
\text{for every $n\in\Z$}.
$$
Clearly the shift operator is an automorphism of $\sC(\cA)$,
and one defines the operator $K^\bullet\mapsto K^\bullet[n]$,
for every $n\in\Z$, as the $n$-th power of the shift operator
(in the automorphism group of $\sC(\cA)$). Then, we can
combine the two previous operators, to define the complex
$$
A[n]:=(A[0])[n]
\qquad
\text{for every $n\in\Z$ and every $A\in\Ob(\cA)$}.
$$
Furthermore, notice that the shift operator restricts
to functors
$$
\sC^+(\cA)\to\sC^+(\cA)
\qquad
\sC^-(\cA)\to\sC^-(\cA)
\qquad
\sC^b(\cA)\to\sC^b(\cA).
$$

\begin{definition}\label{def_homotopies}
Let $\cA$ be an additive category,
$\phi^\bullet,\psi^\bullet:K^\bullet\to L^\bullet$ two morphisms
in $\sC(\cA)$. 

(i)\ \
A {\em (chain) homotopy} from $\phi^\bullet$ to $\psi^\bullet$
is the datum $(s^n:K^n\to L^{n-1}~|~n\in\Z)$ of a system of
morphisms in $\cA$ such that
$$
\phi^n-\psi^n=s^{n+1}\circ d_K^n+d_L^{n-1}\circ s^n
\qquad
\text{for every $n\in\Z$}.
$$

(ii)\ \
We say that $\phi^\bullet$ and $\psi^\bullet$ as in (i) are
{\em chain homotopic} if there is a chain homotopy between
them. It is easily seen that this defines an equivalence
relation $\sim$ on the set $\Hom_{\sC(\cA)}(K^\bullet,L^\bullet)$,
which is preserved by composition of morphisms : if
$\phi^\bullet\sim\psi^\bullet$ and
$\alpha^\bullet:K'{}^\bullet\to K^\bullet$,
$\beta^\bullet:L^\bullet\to L'{}^\bullet$ are any two morphisms,
then $\phi^\bullet\circ\alpha^\bullet\sim\psi^\bullet\circ\alpha^\bullet$
and $\beta^\bullet\circ\phi^\bullet\sim\beta^\bullet\circ\psi^\bullet$.
It follows that, for every interval $I\subset\Z$ there exists
a well defined {\em homotopy category}
$$
\Hot^I(\cA)
$$
whose objects are the same as those of $\sC^I(\cA)$, and whose
morphisms are the homotopy classes of morphisms of complexes,
and a natural functor
\set\begin{equation}\label{eq_cplx-to-hot}
\sC^I(\cA)\to\Hot^I(\cA)
\end{equation}
which is the identity on objects, and the quotient map on
$\Hom$-sets. If $I=\Z$, we write $\Hot(\cA)$ instead of
$\Hot^\Z(\cA)$, and we denote
$$
\rHot_\cA(K^\bullet,L^\bullet)
$$
the set of morphisms $K^\bullet\to L^\bullet$ in $\Hot(\cA)$.
We may also define the subcategories $\Hot^+(\cA)$,
$\Hot^-(\cA)$ and $\Hot^b(\cA)$ as in \eqref{sec_brutal-truncate}.

(iii)\ \
We also say that a morphism $\phi:K^\bullet\to L^\bullet$ is
a {\em homotopy equivalence}, if the class of $\phi^\bullet$
is an isomorphism in $\Hot(\cA)$, {\em i.e.} if there exists
a morphism $\psi^\bullet:L^\bullet\to K^\bullet$ such that
$\psi^\bullet\circ\phi^\bullet\sim\one_{K^\bullet}$ and
$\phi^\bullet\circ\psi^\bullet\sim\one_{L^\bullet}$. We say that
$\phi^\bullet$ is {\em null-homotopic} if it is chain homotopic
to the zero morphism $K^\bullet\to L^\bullet$. We say that
a complex $K^\bullet$ is {\em homotopically trivial}, if the
zero endomorphism $0\cdot\one_{K^\bullet}$ is a homotopy
equivalence.
\end{definition}

\begin{remark}\label{rem_homotopies}
Let $\cA$ and $\cA'$ be any two additive categories.

(i)\ \
If $F:\cA\to\cA'$ is any additive functor, we get
induced functors
$$
\sC(F):\sC(\cA)\to\sC(\cA')
\qquad
\Hot(F):\Hot(\cA)\to\Hot(\cA')
$$
by the rule :
$$
F(K^\bullet)^n:=F(K^n)
\quad\text{and}\quad
d^n_{F(K)}:=F(d^n_K)
\qquad
\text{for every $n\in\Z$ and every $K^\bullet\in\Ob(\sC(\cA))$}.
$$
More generally, we get induced functors $\sC^I(F)$, $\Hot^I(F)$
for any interval $I\subset\Z$, as well as $\sC^+(F)$, $\sC^-(F)$
and $\sC^b(F)$ (and the corresponding homotopy category variants),
in the obvious fashion. Also, notice that the shift operators
descend to a functor on the homotopy category :
$$
\Hot(\cA)\to\Hot(\cA)
\qquad
A^\bullet\mapsto A^\bullet[n]
\qquad
\text{for every $n\in\Z$}
$$
and the latter restrict to endofunctors of the
subcategories $\Hot^+(\cA)$, $\Hot^-(\cA)$ and
$\Hot^b(\cA)$.

(ii)\ \
The category $\sC(\cA)$ is additive. Namely, if
$K^\bullet,L^\bullet$ are any two complexes of $\cA$,
and $\phi^\bullet,\psi^\bullet:K^\bullet\to L^\bullet$
any two morphisms, we define
$$
(\phi+\psi)^i:=\phi^i+\psi^i
\qquad
\text{for every $i\in\Z$}
$$
and it is clear that this rule yields a natural abelian
group structure on $\Hom_\cA(K^\bullet,L^\bullet)$, such
that composition of morphisms is a bilinear operation.
The zero object is the (unique) complex $0$ which has
the zero object of $\cA$ in each degree, and biproducts
in $\sC(\cA)$ admit natural representatives, defined by
the rule :
$$
(A\oplus B)^i:=A^i\oplus B^i
\qquad
\text{for every $A^\bullet,B^\bullet\in\Ob(\sC(\cA))$ and
every $i\in\Z$}.
$$
Moreover, it is clear that if
$\phi_1^\bullet,\phi_2^\bullet:K^\bullet\to L^\bullet$ are
any two morphisms in $\sC(\cA)$ such that
$\phi^\bullet_1\sim\phi^\bullet_2$, we have
$\phi^\bullet_1+\psi^\bullet\sim\phi^\bullet_2+\psi^\bullet$
for every other morphism $\psi^\bullet:K^\bullet\to L^\bullet$;
hence there is a unique abelian group structure on
$\rHot_\cA(K^\bullet,L^\bullet)$ for every
$K^\bullet,L^\bullet\in\Ob(\sC(\cA))$, that makes
$\Hot(\cA)$ into an additive category, such that
\eqref{eq_cplx-to-hot} is an additive functor.

Furthermore, if $\cA$ is abelian, the same holds
for $\sC(\cA)$. Indeed, we can take
$$
(\Ker\,\phi)^i:=\Ker\,\phi^i
\qquad
(\Coker\,\phi)^i:=\Coker\,\phi^i
\qquad
\text{for every morphism $\phi^\bullet$ of $\sC(\cA)$}
$$
and it is clear that the resulting morphism
$\beta_{\phi^\bullet}$ as in \eqref{eq_ker-coker} is
an isomorphism. It is then also clear that, for any
additive functor $F$ as in (i), the induced functors
$\sC(F)$ and $\Hot(F)$ are additive, and the same holds
as well for the shift operators on $\sC(\cA)$ and
$\Hot(\cA)$.

All these considerations extend more generally to
the full subcategories $\sC^I(\cA)$ (for any interval
$I\subset\Z$), $\sC^-(\cA)$, $\sC^+(\cA)$ and
$\sC^b(\cA)$, as well as their homotopy category
variants : the detailed verifications shall
be left to the reader.

(iii)\ \
Taking into account remark \ref{rem_Add.Fun}(i),
we get a natural identification
$$
\sC(\cA)^o\isom\sC(\cA^o)
\qquad
K^\bullet\mapsto(K^o)^\bullet
$$
where
$$
(K^o)^i:=(K^{-i})^o
\qquad\text{and}\qquad
d^i_{K^o}:=(-1)^{i+1}\cdot(d_K^{-i-1})^o
$$
for every $i\in\Z$ and every $K^\bullet\in\Ob(\sC(\cA))$.
If $\phi^\bullet:K^\bullet\to L^\bullet$ is any morphism in
$\sC(\cA)$, then $(\phi^o)^i:=\phi^{-i}$ for every $i\in\Z$.
Also, if $s^\bullet$ is a homotopy between morphisms
$f^\bullet,g^\bullet:K^\bullet\to L^\bullet$, then the system
$$
(s^o)^\bullet:=((-1)^i\cdot(s^{-i})^o~|~i\in\Z)
$$
is a homotopy between
$(f^o)^\bullet,(g^o)^\bullet:(L^o)^\bullet\to(K^o)^\bullet$,
whence an isomorphism of categories :
$$
\Hot(\cA)^o\isom\Hot(\cA^o).
$$
One may wonder about the reason for inserting apparently
artificial signs in the above formulas; the point is that
this sign convention makes the foregoing natural
identifications behave slightly better in connection with
additive contravariant functor (and their extensions to
complexes : see {\em e.g.} the construction of the
$\Hom^\bullet$ functor in example \ref{ex_dual-is-opposite}).
On the other hand, there is a drawback as well : indeed,
notice that by applying twice our natural isomorphism,
we do {\em not} get the identity automorphism of
$\sC(\cA)^{oo}=\sC(\cA)$; rather
$$
d_{K^{oo}}^\bullet=-d_K^\bullet
\qquad\text{and}\qquad
(s^{oo})^\bullet=s^\bullet
$$
for every object $K^\bullet$ and every homotopy $s^\bullet$
of $\sC(\cA)$ (details left to the reader).

(iv)\ \
The indexing notation that makes use of superscript to
denote the degrees in a complex, is known traditionally
as {\em cohomological degree notation}. Sometimes it is
more natural to switch to the {\em homological degree
notation}, that makes use of subscript indexing; namely,
one associates with any cochain complex $K^\bullet$, the
{\em chain complex} $K_\bullet$ given by the rule :
$$
K_n:=K^{-n}
\qquad\text{and}\qquad
d_n:=d^{-n}:K_n\to K_{n-1}
\qquad
\text{for every $n\in\Z$}.
$$

(v)\ \
If $\cA$ is complete (resp. cocomplete) then the same
holds for $\sC^I(\cA)$, for any interval $I\subset\Z$.
Indeed, we may regard these categories as full subcategories
of the category $\Fun(\Z,\cA)$, where the ordered set
$\Z$ (endowed with its standard ordering) is viewed as
a category as explained in example \ref{ex_universe}(iii),
and therefore the assertion follows from corollary
\ref{lem_presheaf-in-a-cat}(ii,iv), which shows more
precisely that the limits and colimits in $\sC^I(\cA)$
are formed {\em degree-wise}.

(vi)\ \
If $\cA$ is complete (resp. cocomplete), then all
products (resp. coproducts) of $\Hot^I(\cA)$ are
representable, for every interval $I\subset\Z$.
Indeed, let $K^\bullet_\bullet:=(K^\bullet_j~|~j\in J)$
be any family of objects of $\Hot^I(\cA)$ indexed by
a small set $J$; we claim that the product $K^\bullet$ in
the category $\sC^I(\cA)$ of the family $K^\bullet_\bullet$
also represents the product of the same family in
$\Hot^I(\cA)$. Indeed, let $M^\bullet$ be any other object
of $\Hot^I(\cA)$, and $\pi^\bullet_j:K^\bullet\to K^\bullet_j$
the canonical projection, for every $j\in J$; from any
morphism $\phi^\bullet:M^\bullet\to K^\bullet$ we obtain a
system
$(\phi^\bullet_j:=
\pi^\bullet_j\circ\phi^\bullet:M^\bullet\to K^\bullet_j~|~j\in J)$
of morphisms in $\Hot^I(\cA)$. Conversely, given such a
system of morphisms, for every $j\in J$ pick arbitrarily
a morphism of complexes
$\tilde\phi^\bullet_j:M^\bullet\to K^\bullet_j$ in the class of
$\phi^\bullet_j$; the system $(\tilde\phi^\bullet_j~|~j\in J)$
yields by (v) a unique morphism of complexes
$\tilde\phi^\bullet:M^\bullet\to K^\bullet$ and we let
$\phi^\bullet$ be the homotopy class of $\tilde\phi^\bullet$.
We need to show that $\phi^\bullet$ does not depend on the
choices of representatives $\tilde\phi^\bullet_j$; to this
aim, we come down to checking that if each $\tilde\phi^\bullet_j$
is null-homotopic, then the same holds for $\tilde\phi^\bullet$.
Hence, pick for every $j\in J$ a homotopy $s^\bullet_j$
from $\tilde\phi^\bullet_j$ to the zero morphism. Since
the products are computed degree-wise in $\sC^I(\cA)$,
the system $(s^n_j:M^n\to K^{n-1}_j~|~j\in J)$ corresponds,
for every $n\in I$, to a unique morphism $s^n:M^n\to K^{n-1}$
of $\cA$, and it is easily seen that the system $(s^n~|~n\in I)$
yields a well defined homotopy from $\tilde\phi^\bullet$ to the
zero morphism $M^\bullet\to K^\bullet$. By applying the same
argument to the opposite category $\cA^o$ and invoking (iii)
we get the assertion for coproducts.
\end{remark}

\begin{example}\label{ex_dual-is-opposite}
Denote by $\Z\Mod_\mathrm{fft}$ the additive category of
free abelian groups of finite rank. We have a natural
isomorphism of categories :
$$
\Z\Mod_\mathrm{fft}^o\isom\Z\Mod_\mathrm{fft}
\quad :\quad
P\mapsto P^\vee:=\Hom_\Z(P,\Z)
$$
that sends any $\Z$-linear map $f:P\to Q$ to its
transpose $f^\vee:Q^\vee\to P^\vee$. In view of remark
\ref{rem_homotopies}(iii) there follow natural
isomorphisms of categories
$$
\sC(\Z\Mod_\mathrm{fft})^o\isom\sC(\Z\Mod_\mathrm{fft})
\qquad
\Hot(\Z\Mod_\mathrm{fft})^o\isom\Hot(\Z\Mod_\mathrm{fft})
\qquad
K^\bullet\mapsto K^{\vee\bullet}.
$$
\end{example}

\sset\subsubsection{}
A {\em double complex} of $\cA$ is an object of
$$
\sC_2(\cA):=\sC(\sC(\cA))
$$
(and likewise for a morphism of double complexes);
{\em i.e.} a double complex is a triple
$$
K^{\bullet\bullet}:=
(K^{\bullet\bullet},d^{\bullet\bullet}_h,d^{\bullet\bullet}_v)
$$
consisting of a system $(K^{pq}~|~p,q\in\Z)$ of objects of
$\cA$, and morphisms $d^{pq}_h$, $d^{pq}_v$ called respectively
the {\em horizontal} and {\em vertical} differentials, fitting
into a commutative diagram
$$
{\diagram
K^{pq} \ar[r]^-{d^{pq}_h} \ar[d]_{d^{pq}_v} & K^{p+1,q} \ar[d]^{d_v^{p+1,q}} \\
K^{p,q+1} \ar[r]^-{d^{pq}_h} & K^{p+1,q+1}
\enddiagram}
\qquad
\text{for every $p,q\in\Z$}
$$
and such that
$$
d^{p+1,q}_h\circ d^{pq}_h=0
\qquad
d^{p,q+1}_v\circ d^{pq}_v=0
\qquad
\text{for every $p,q\in\Z$}.
$$

\sset\subsubsection{}\label{subsec_flip-diag}
There are natural functors
$$
\sC_2(\cA)\to\sC_2(\cA)
\ :\
K^{\bullet\bullet}\mapsto\mathsf{fl}_\cA(K^{\bullet\bullet})
\qquad\qquad
\sC_2(\cA)\to\sC(\cA)
\ :\
K^{\bullet\bullet}\mapsto(K^{\bullet\bullet})^\Delta
$$
where :
\begin{itemize}
\item
The {\em flip} $\mathsf{fl}_\cA(K^{\bullet\bullet})$ of
$K^{\bullet\bullet}$ is the double complex $F^{\bullet\bullet}$
such that $F^{pq}:=K^{qp}$ for every $p,q\in\Z$, with differentials
deduced from those of $K^{\bullet\bullet}$, in the obvious way.
\item
The {\em diagonal} $(K^{\bullet\bullet})^\Delta$ is
the complex $D^\bullet$ such that $D^p:=K^{pp}$ for every
$p\in\Z$ and with differentials
$$
d^{p+1,q}_v\circ d^{pq}_h:D^p\to D^{p+1}
\qquad
\text{for every $p\in\Z$}.
$$
\end{itemize}
Suppose that all coproducts (resp. all products) are representable
in $\cA$. Then there are two other natural functors
$$
\Tot_\cA^\oplus:\sC_2(\cA)\to\sC(\cA)
\qquad
\text{(\ resp. $\Tot_\cA^\Pi:\sC_2(\cA)\to\sC(\cA)$ \ )}
$$
defined as follows. The {\em total complex}
$\Tot^\oplus_\cA(K^{\bullet\bullet})$ (resp.
$\Tot^\Pi_\cA(K^{\bullet\bullet})$) is the complex $T^\bullet$
with
$$
T^n:=\bigoplus_{p+q=n}K^{pq}
\qquad
\text{( resp. $T^n:=\prod_{p+q=n}K^{pq}$ )}
$$
and with differentials $T^n\to T^{n+1}$ given by the sum
(resp. the product) of the morphisms
$$
d^{pq}_h+(-1)^p\cdot d^{pq}_v:K^{pq}\to K^{p,q+1}\oplus K^{p+1,q}
\qquad
\text{for all $p,q\in\Z$ such that $p+q=n$}.
$$
We usually drop the subscript $\cA$, unless the omission
would be a source of ambiguities, and we often omit as
well the superscript $\oplus$ when dealing with the total
complex functor; to avoid confusion, we stipulate that
{\em the notation $\Tot$ shall always refer to the
functor $\Tot^\oplus$}, so if we need to use the other
total complex functor, we shall denote it explicitly by
$\Tot^\Pi$. Notice that we have natural isomorphisms
$$
\Tot(K^{\bullet\bullet})\isom\Tot(\mathsf{fl}_\cA(K^{\bullet\bullet}))
\qquad
\text{(resp.\ 
$\Tot^{\Pi}(K^{\bullet\bullet})\isom
\Tot^{\Pi}(\mathsf{fl}_\cA(K^{\bullet\bullet}))$
\ )}
$$
given, in each degree $n\in\Z$, by the direct sum (resp. the
direct product) of the morphisms $(-1)^{pq}\cdot\one_{K^{pq}}$,
for every $p,q\in\Z$ such that $p+q=n$.

\sset\subsubsection{}\label{subsec_triple-complex}
The category of {\em triple complexes} of $\cA$ can be
realized in two equivalent ways : namely, we have natural
identifications :
\set\begin{equation}\label{eq_triple-cplx}
\sC_2(\sC(\cA))\isom\sC(\sC_2(\cA))
\end{equation}
and we denote either of these categories by $\sC_3(\cA)$.
In other words, a triple complex is a system
$$
(K^{ijk},d_1^{ijk},d_2^{ijk},d_3^{ijk}~|~i,j,k\in\Z)
$$
such that, for every fixed $p\in\Z$, the subsystems
$$
(K^{p,\bullet\bullet},d_2^{p,\bullet\bullet},d_3^{p,\bullet\bullet})
\qquad
(K^{\bullet,p,\bullet},d_1^{\bullet,p,\bullet},d_3^{\bullet,p,\bullet})
\qquad
(K^{\bullet,\bullet,p},d_1^{\bullet,\bullet,p},d_2^{\bullet,\bullet,p})
$$
are double complexes of $\cA$. Corresponding to the
two interpretations of $\sC_3(\cA)$, we have two different
way of forming total complexes out of triple complexes;
namely, if $\cA$ is cocomplete we have the functors
$$
\Tot_\cA\circ\Tot_{\sC(\cA)},\Tot_\cA\circ\sC(\Tot_\cA):
\sC_3(\cA)\to\sC(\cA)
$$
and likewise for the $\Tot^\Pi$ variant, if $\cA$ is
complete. However, a simple computation shows that these
two functors agree under the natural identification
\eqref{eq_triple-cplx} : namely, both functors
yield a complex whose differential is the direct
sum (resp. direct product) of the morphisms
$$
d_1^{ijk}+(-1)^i\cdot d_2^{ijk}+(-1)^{i+j}\cdot d_3^{ijk}:
K^{ijk}\to K^{i+1,j,k}\oplus K^{i,j+1,k}\oplus K^{i,j,k+1}.
$$
We shall therefore denote these functors indifferently
by $\Tot$ (resp. by $\Tot^\Pi$), just as for double complexes.

\sset\subsubsection{}\label{subsec_biadditive}
Let $\cA,\cA',\cA''$ be three additive categories, and
$$
B:\cA\times\cA'\to\cA''
$$
a {\em biadditive} functor, {\em i.e.} such that, for
every $A\in\Ob(\cA)$ and every $A'\in\Ob(\cA')$, both
functors
$$
B(A,-):\cA'\to\cA''
\qquad
B(-,A'):\cA'\to\cA''
$$
are additive. The functor $B$ induces a functor
$$
B^{\bullet\bullet}:\sC(\cA)\times\sC(\cA')\to\sC_2(A'')
$$
by the rule
$$
B^{ij}(K^\bullet,L^\bullet):=B(K^i,L^j)
\qquad
\text{for every $K^\bullet\in\Ob(\sC(\cA))$ and
$L^\bullet\in\Ob(\sC(\cA'))$}.
$$
The differentials of $B^{\bullet\bullet}(K^\bullet,L^\bullet)$
in degree $(p,q)$ are respectively $B(d^p_K,\one_{L^q})$ and
$B(\one_{K^p},d^q_L)$. Suppose that $\cA''$ is cocomplete
(resp. complete); then, by composing with the functor
$\Tot^\oplus$ (resp. with $\Tot^\Pi$), we deduce a biadditive
functor
$$
B^\bullet_\oplus:\sC(\cA)\times\sC(\cA')\to\sC(\cA'')
\qquad
\text{(resp.\ $B^\bullet_\Pi:\sC(\cA)\times\sC(\cA')\to\sC(\cA'')$)}.
$$
For a general additive category $\cA''$, the foregoing functors
are still well defined at least on bounded above or bounded
below complexes, {\em i.e.} we have functors
$$
B^\bullet_{-}:\sC^-(\cA)\times\sC^-(\cA')\to\sC^-(\cA'')
\qquad
B^\bullet_{+}:\sC^+(\cA)\times\sC^+(\cA')\to\sC^+(\cA'').
$$

\begin{example}\label{ex_Hom-complex-hots}
(i)\ \
If $\cA$ is any additive category, the construction of
\eqref{subsec_biadditive} extends the biadditive functor
$$
\Hom_\cA(-,-):\cA^o\times\cA\to\Z\Mod
$$
to a functor
$\Hom^\bullet_{\cA,\Pi}:\sC(\cA)\times\sC(\cA^o)\to\sC(\Z\Mod)$
(notice the swapping of the two factors, that is required
in order to conform with established sign conventions)
and combining with the natural identification of remark
\ref{rem_homotopies}(iii), we obtain a biadditive functor
$$
\Hom^\bullet_\cA(-,-):\sC(\cA)\times\sC(\cA)^o\to\sC(\Z\Mod)
\quad :\quad
(L^\bullet,K^\bullet)\to\Hom_{\cA,\Pi}^\bullet((K^o)^\bullet,L^\bullet).
$$
We shall denote by $d^\bullet_{K,L}$ the differential of
the complex $\Hom^\bullet_\cA(K^\bullet,L^\bullet)$, for every
$K^\bullet,L^\bullet\in\Ob(\sC(\cA))$. Let $n\in\Z$ be any
integer; with this notation, a simple inspection shows that :
\begin{itemize}
\item
$d^n_{K,L}$ is the product, for every $j\in\Z$, of the maps
$$
\Hom_\cA(K^{j-n},L^j)\to
\Hom_\cA(K^{j-n},L^{j+1})\oplus\Hom_\cA(K^{j-n-1},L^j)
$$
given by the rule :
$(\phi:K^{j-n}\to L^j)\mapsto
d^j_L\circ\phi-(-1)^n\cdot\phi\circ d^{j-n-1}_K$.
\item
$\Hom_{\sC(\cA)}(K^\bullet,L^\bullet[n])=\Ker\,d^n_{K,L}$.
\item
Let $\phi^\bullet,\psi^\bullet:K^\bullet\to L^\bullet[n]$ be any
two morphisms; then the set of homotopies from
$\phi^\bullet$ to $\psi^\bullet$ is naturally identified
with the subset
$$
\{s^\bullet\in\Hom_\cA^{-1-n}(K^\bullet,L^\bullet)
~|~d^{-1-n}_{K,L}(s^\bullet)=\psi^\bullet-\phi^\bullet\}.
$$
\item
Consequently, we have a natural identification :
$$
\rHot_\cA(K^\bullet,L^\bullet[n])\isom
H^n\Hom^\bullet_\cA(K^\bullet,L^\bullet).
$$
\end{itemize}

(ii)\ \
The constructions of (i) can be used to endow $\sC(\cA)$
with a natural $2$-category structure, whose $2$-cells are
given by homotopies of complexes. Indeed, let us write
$$
s^\bullet:\phi^\bullet\Rightarrow\psi^\bullet
$$
if $s^\bullet:=(s^n~|~n\in\N)$ is a homotopy from $\phi^\bullet$
to $\psi^\bullet$. Then, if $\lambda^\bullet:K^\bullet\to L^\bullet$
is another morphism, and
$t^\bullet:\psi^\bullet\Rightarrow\lambda^\bullet$ another homotopy,
we define a composition law by setting
$$
s^\bullet\odot t^\bullet:=(s^n+t^n~|~n\in\N)
$$
and it is immediate that $s^\bullet\odot t^\bullet$ is a homotopy
$\phi^\bullet\Rightarrow\lambda^\bullet$. Moreover, if
$\beta^\bullet,\gamma^\bullet:L^\bullet\to P^\bullet$ are two other
morphisms in $\sC(\cA)$, and $u:\beta^\bullet\Rightarrow\gamma^\bullet$
another homotopy, we have a Godement composition law by the rule
$$
u^\bullet*s^\bullet:=(\beta^{n+1}\circ s^n+u^n\circ\psi^n~|~n\in\N)
:\beta^\bullet\circ\phi^\bullet\Rightarrow\gamma^\bullet\circ\psi^\bullet.
$$
The associativity of the laws $*$ and $\odot$ thus defined
are easily checked by direct computation.

Now, suppose that $\delta^\bullet:L^\bullet\to P^\bullet$ is yet another
morphism, and $v:\gamma\Rightarrow\delta$ another homotopy.
A direct calculation yields the identity
$$
(u^\bullet*s^\bullet)\odot(v^\bullet*t^\bullet)=
(u^\bullet\odot v^\bullet)*(s^\bullet\odot t^\bullet)+c^\bullet
$$
where $c^n:=d^{n-2}_P\circ u^{n-1}\circ t^n-u^n\circ t^{n+1}\circ d^n_K$
for every $n\in\Z$, which can be rewritten as
$$
c^\bullet=d^{-2}_{K,L}((u\circ t)^\bullet)
\qquad
\text{where $(u\circ t)^n:=u^{n-1}\circ t^n$ for every $n\in\Z$}.
$$
Summing up, we conclude that $\sC(\cA)$ carries a $2$-category
structure, whose $2$-cells $\phi\Rightarrow\psi$ (for any two
$1$-cells $\phi,\psi:K^\bullet\to L^\bullet$) are the classes
$$
\bar s\in\Hom^{-1}_\cA(K^\bullet,L^\bullet)/\Img(d^{-2}_{K,L})
\qquad
\text{such that $d^{-1}_{K,L}(\bar s)=\psi^\bullet-\phi^\bullet$}.
$$

(iii)\ \
Moreover, if $F$ is any additive functor as in remark
\ref{rem_homotopies}(i), the induced functor $\sC(F)$
extends to a pseudo-functor for the $2$-category
structures given by (ii); indeed, if
$s^\bullet:K^\bullet\Rightarrow L^\bullet$ is a homotopy,
obviously the system $(Fs^n~|~n\in\N)$ is a homotopy
$Fs^\bullet:F(K^\bullet)\Rightarrow F(L^\bullet)$.
\end{example}

\begin{definition}
Let $\cA$ be an abelian category, $K^\bullet$ any object
of $\sC(\cA)$, and  $i\in\Z$ any integer.
\begin{enumerate}
\item
The {\em cohomology of $K^\bullet$ in degree $i$} is the object
of $\cA$ :
$$
H^iK^\bullet:=\Ker\,d^i/\Img\,d^{i-1}.
$$
\item
We say that $K^\bullet$ is {\em acyclic in degree $i$}
(resp. {\em acyclic}) if $H^iK^\bullet=0$ (resp. if
$H^jK^\bullet=0$ for every $j\in\Z$).
\end{enumerate}
\end{definition}

\begin{remark}\label{rem_abel-homotopies}
(i)\ \
Clearly, for every $i\in\Z$, the rule
$K^\bullet\mapsto H^iK^\bullet$ extends to a functor
$$
H^i:\sC(\cA)\to\cA.
$$

(ii)\ \
Moreover, if $\phi^\bullet,\psi^\bullet:K^\bullet\to L^\bullet$
are chain homotopic morphisms in $\sC(\cA)$, then it is easily
seen that, for every $i\in\Z$, the induced morphisms in cohomology
$$
H^i\phi^\bullet,H^i\psi^\bullet:H^iK^\bullet\to H^iL^\bullet
$$
coincide. Hence, the cohomology functor $H^i$ factors (uniquely)
through a functor
$$
H^i:\Hot(\cA)\to\cA
\qquad
\text{for every $i\in\Z$}.
$$

(iii)\ \
Notice as well that the cohomology is a {\em self-dual}
operation : namely, we have an essentially commutative
diagram of functors
$$
{\diagram \sC(\cA)^o \ar[rr]^-\sim \ar[rd]_{(H^i)^o} & &
\sC(\cA^o) \ar[ld]^{H^{-i}} \\
& \cA^o
\enddiagram}
\qquad
\text{for every $i\in\Z$}
$$
whose top horizontal arrow is the isomorphism of
remark \ref{rem_homotopies}(iii). Indeed, if $K^\bullet$
is any complex of $\cA$, the morphism $d^{i-1}_K$
factors uniquely through a morphism
$e^{i-1}_K:K^{i-1}\to\Ker\,d^i_K$, and we have a natural
isomorphism $H^iK^\bullet\isom\Coker\,e^{i-1}_K$. On the
other hand, the morphism $d^i_K$ factors uniquely
through a morphism $f^i_K:\Coker\,d^{i-1}_K\to K^i$
and we have also the natural isomorphism
$H^iK^\bullet\isom\Ker\,f^i_K$, so the assertion
follows from remark \ref{rem_Add.Fun}(i).

(iv)\ \
Just like in remark \ref{rem_homotopies}(iv), it
is sometimes convenient to switch to a subscript
notation; thus, for any chain complex $K_\bullet$
one sets
$$
H_nK_\bullet:=H^{-n}K_n
\qquad
\text{for every $n\in\Z$}
$$
and calls this object of $\cA$ the {\em homology of $K_\bullet$}
in degree $n$.
\end{remark}

\begin{example}\label{ex_monoidal}
(i)\ \ 
Suppose that $(\cA,\otimes,\Phi,\Psi)$ is a tensor abelian
category, and set
$$
B(A,A'):=A\otimes A'
\qquad
\text{for every $A,A'\in\Ob(\cA)$}.
$$
Following \eqref{subsec_biadditive}, we get a corresponding
functor $B^{\bullet\bullet}_\oplus$, for which we use the
notation
$$
-\boxtimes-:\sC(\cA)\times\sC(\cA)\to\sC_2(\cA).
$$
If all coproducts are representables in $\cA$, we get
as well the functor $B^\bullet_\oplus$, which we denote :
$$
K^\bullet\otimes L^\bullet:=\Tot(K^\bullet\boxtimes L^\bullet).
$$
If $\cA$ is not cocomplete, $K^\bullet\otimes L^\bullet$
is still well defined, provided $K^\bullet$ and $L^\bullet$
are both bounded below or if they are both bounded above.
The commutativity constraints for $(\cA,\otimes)$ yield
natural isomorphisms $K^\bullet\boxtimes L^\bullet\isom
\mathsf{fl}_\cA(L^\bullet\boxtimes K^\bullet)$, as well as
\set\begin{equation}\label{eq_commut-constraint}
\Psi^\bullet_{K,L}:
K^\bullet\otimes L^\bullet\isom L^\bullet\otimes K^\bullet.
\end{equation}
Namely, one takes the direct sum of the maps
$(-1)^{pq}\cdot\Psi_{K^p,L^q}$, for every $p,q\in\Z$.

(ii)\ \
Likewise, if $P^\bullet$ is another complex of $\cA$, the
associativity constraints of $\cA$ assemble to an isomorphism
of triple complexes
$$
(\Phi_{K^i,L^j,P^k}~|~i,j,k):
(K^\bullet\boxtimes L^\bullet)\boxtimes P^\bullet\isom
K^\bullet\boxtimes(L^\bullet\boxtimes P^\bullet)
$$
whence, taking into account the discussion of
\eqref{subsec_triple-complex}, a natural isomorphism
in $\sC(\cA)$ :
$$
\Phi^\bullet_{K,L,P}:
K^\bullet\otimes(L^\bullet\otimes P^\bullet)\isom
(K^\bullet\otimes L^\bullet)\otimes P^\bullet.
$$
With these natural isomorphisms, $\sC(\cA)$ is then
a tensor abelian category as well.

(iii)\ \
In the situation of (i), notice that the natural morphism
$K^i\otimes L^j\to(K^\bullet\otimes L^\bullet)^{i+j}$ induces
morphisms
$$
\begin{aligned}
&\Ker\,(d^i_K)\otimes\Ker\,(d^j_L)\to\Ker\,(d^{i+j}_{K\otimes L}) \\
&(\Ker\,(d^i_K)\otimes\Img\,(d^{j-1}_L))\oplus
(\Img\,(d^{i-1}_K)\otimes\Ker\,(d^j_L))\to\Img\,(d^{i+j}_{K\otimes L})
\end{aligned}
$$
so, the induced map
$\Ker\,(d^i_K)\otimes\Ker\,(d^j_L)\to H^{i+j}(K^\bullet\otimes L^\bullet)$
factors through a natural pairing :
$$
H^i(K^\bullet)\otimes H^j(L^\bullet)\to H^{i+j}(K^\bullet\otimes L^\bullet)
\qquad
\text{for every $i,j\in\Z$}.
$$

(iv)\ \
Moreover, if $P^\bullet$ is a third complex, in view of
(ii) we get a commutative diagram
$$
\xymatrix{
H^iK^\bullet\otimes(H^jL^\bullet\otimes H^kP^\bullet)
\ar[d]_\gamma \ar[r]^-\alpha & \ar[r]^-\beta
(H^iK^\bullet\otimes H^jL^\bullet)\otimes H^kP^\bullet &
H^{i+j}(K^\bullet\otimes L^\bullet)\otimes H^kP^\bullet \ar[d]^\delta \\
H^iK^\bullet\otimes H^{j+k}(L^\bullet\otimes P^\bullet) \ar[r]^\lambda
& H^{i+j+k}(K^\bullet\otimes(L^\bullet\otimes P^\bullet)) \ar[r]^-\sigma
& H^{i+j+k}((K^\bullet\otimes L^\bullet)\otimes P^\bullet)
}$$
where $\alpha$ is the associativity constraint, $\beta$,
$\gamma$, $\delta$ and $\lambda$ are given by the above
pairing, and $\sigma$ is deduced from $\Phi^\bullet_{K,L,P}$.

(v)\ \
If $\cA$ also admits an internal $\Hom$ functor, we may define
as well a functor
$$
\cHom^{\bullet\bullet}:
\sC(\cA)\times\sC(\cA)^o\to\sC_2(\cA)
$$
following the trace of example \ref{ex_Hom-complex-hots}(i);
namely, we set
$$
\cHom^{p,q}(K^\bullet,L^\bullet):=\cHom(K^{-p},L^q)
\qquad
\text{for every $p,q\in\Z$}
$$
with differentials
$$
d_h^{p,q}:=\cHom(\one_{K^{-p}},d_L^q)
\qquad
d_v^{p,q}:=(-1)^{p+1}\cdot\cHom(d_K^{-p-1}\one_{L^q})
\qquad
\text{for every $p,q\in\Z$}.
$$
If all products are representable in $\cA$, we may then define
$$
\cHom^\bullet(K^\bullet,L^\bullet):=
\Tot^\Pi\cHom^{\bullet\bullet}(K^\bullet,L^\bullet).
$$

(vi)\ \
Suppose that $\cA$ is both complete and cocomplete,
and $P^\bullet$ is any other complex; to ease notation, set
$H^\bullet:=\cHom^\bullet(K^\bullet,L^\bullet)$ and
$N^\bullet:=P^\bullet\otimes K^\bullet$. We have natural
isomorphisms
$$
\begin{aligned}
\Hom_{\sC(\cA)}(P^\bullet,H^\bullet)\isom\, &
\xymatrix{\Equal\Bigl(\Hom^0_\cA(P^\bullet,H^\bullet)
\ar@<.5ex>[r]^-{d_v} \ar@<-.5ex>[r]_-{d_h} &
\Hom^1_\cA(P^\bullet,H^\bullet)\Bigr)}
\end{aligned}
$$
where
$$
d_h:=\prod_{n\in\Z}\Hom_\cA(P^n,d_H^n)
\qquad\text{and}\qquad
d_v:=\prod_{n\in\Z}\Hom_\cA(d_P^n,H^{n+1})
$$
(see example \ref{ex_Hom-complex-hots}(i)). However, notice that
$$
\begin{aligned}
\Hom^a_\cA(P^\bullet,H^\bullet)=\, & 
\prod_{n\in\Z}\Hom_\cA\Bigl(P^n,\prod_{p+q=n+a}\cHom(K^{-p},L^q)\Bigr) \\
=\, & \prod_{n\in\Z}\ \prod_{p+q=n+a}\Hom_\cA(P^n,\cHom(K^{-p},L^q)) \\
=\, & \prod_{n\in\Z}\ \prod_{p+q=n+a}\Hom_\cA(P^n\otimes K^{-p},L^q) \\
=\, & \prod_{q\in\Z}\Hom_\cA\Bigl(\bigoplus_{n-p=q-a}P^n\otimes K^{-p},L^q\Bigr) 
\end{aligned}
$$
for every $n,a\in\Z$, whence
$$
\Hom_{\sC(\cA)}(P^\bullet,H^\bullet)\isom
\xymatrix{\Equal\Bigl(\Hom^0_\cA(N^\bullet,L^\bullet)
\ar@<.5ex>[r]^-{d'_v} \ar@<-.5ex>[r]_-{d'_h} &
\Hom^1_\cA(N^\bullet,L^\bullet)\Bigr)}
$$
where
$$
d'_h:=\prod_{n\in\Z}\Hom_\cA(d^q_N,L^q)
\qquad\text{and}\qquad
d'_v:=\prod_{n\in\Z}\Hom_\cA(N^q,d^q_L)
$$
so finally :
$$
\Hom_{\sC(\cA)}(P^\bullet,\cHom^\bullet(K^\bullet,L^\bullet))\isom
\Hom_{\sC(\cA)}(P^\bullet\otimes K^\bullet,L^\bullet)
$$
which says that $\cHom^\bullet$ is an internal $\Hom$ functor
for $\sC(\cA)$.

(vii)\ \
In the situation of (vi), set
$$
\cHom_{\sC(\cA)}(K^\bullet,L^\bullet):=\Ker\,(d^0:H^0\to H^1)
$$
and take $P^\bullet:=Z[0]$, where $Z$ is any object of $\cA$;
it is easily seen that the natural map
$$
\Hom_\cA(Z,\cHom_{\sC(\cA)}(K^\bullet,L^\bullet))\to
\Hom_{\sC(\cA)}(P^\bullet,H^\bullet)\to
\Hom_{\sC(\cA)}(Z[0]\otimes K^\bullet,L^\bullet)
$$
is an isomorphism : details left to the reader.
\end{example}

\sset\subsubsection{}\label{subsec_cplx-mixed-tens}
Suppose that $\cA$ is an additive category with small
$\Hom$-sets, and let $(G^\bullet,d^\bullet_G)$ be any complex
of free abelian groups of finite rank; with the notation
of \eqref{subsec_mixed-tensors}, we obtain an object
$(G^\bullet_\cA,d^\bullet_\cA)$ of $\sC(\cA^\dagger)$;
on the other hand, if $(K^\bullet,d^\bullet_K)$ is any
object of $\sC(\cA)$, we may also consider the object
$h^\dagger_K:=(h^\dagger_{K^n},h^\dagger_{d^n}~|~n\in\N)$ of
$\sC(\cA^\dagger)$, and since $\cA^\dagger$ is an abelian
tensor category, we may form the tensor product
$$
G^\bullet\boxtimes_\Z K^\bullet:=\cHom^{\bullet\bullet}(G^\bullet_\cA,h^\dagger_K)
$$
according to example \ref{ex_monoidal}(v). Arguing as in
\eqref{subsec_mixed-tensors}, we see that this object
of $\sC_2(\cA^\dagger)$ is isomorphic to an object of
$\sC_2(\cA)$, and after choosing representing objects,
we get a functor
$$
\sC(\Z\Mod_\mathrm{fft})\times\sC(\cA)\to\sC_2(\cA)
\qquad
(G^\bullet,K^\bullet)\mapsto G^\bullet\boxtimes_\Z K^\bullet
$$
which is additive in both arguments. If all direct sums
are representable in $\cA$, or else, if $G^\bullet$ is a
bounded complex, it is then natural to define
$$
G^\bullet\otimes_\Z K^\bullet:=
\Tot\,G^\bullet\boxtimes_\Z K^\bullet
\qquad
\text{for every $G^\bullet$ and $K^\bullet$ as above}.
$$
Likewise, we set
$K^\bullet\boxtimes_\Z G^\bullet:=
\mathsf{fl}_\cA(G^\bullet\boxtimes_\Z K^\bullet)$ and
$K^\bullet\otimes_\Z G^\bullet:=
\Tot(K^\bullet\boxtimes_\Z G^\bullet)$.

\begin{remark}\label{rem_chain-homotopies}
(i)\ \
With the notation of \eqref{subsec_cplx-mixed-tens},
notice the natural isomorphism
$$
K^\bullet[1]\isom\Z[1]\otimes_\Z K^\bullet
\qquad
\text{for every $K^\bullet\in\Ob(\sC(\cA))$}
$$
which explains the sign convention in the definition of the
shift operator in \eqref{subsec_shift}. 

(ii)\ \
Moreover, denote by $\K\La 1\Ra^\bullet\in\Ob(\sC^{[-1,0]}(\Z\Mod))$
the object such that
$$
\K\La 1\Ra^{-1}:=\Z
\qquad
\K\La 1\Ra^0:=\Z\oplus\Z
$$
and with differential $d^{-1}$ given by the rule : $n\mapsto(n,-n)$
for every $n\in\Z$. Let $e_0:=(0,1)$ and $e_1:=(1,0)$ be the
canonical basis of $\K\La 1\Ra^0$; we have two morphisms
$$
\iota^\bullet_i:\Z[0]\to\K\La 1\Ra^\bullet
\qquad
\text{for $i=0,1$}
$$
given, in degree zero, by the rule : $n\mapsto n\cdot e_i$ for
every $n\in\N$. For any pair of morphisms
$\phi^\bullet,\psi^\bullet:L^\bullet\to M^\bullet$ in $\sC(\cA)$,
and any homotopy $s^\bullet:=(s^n~|~n\in\Z)$ from $\phi^\bullet$
to $\psi^\bullet$, we obtain a morphism of complexes
$$
\sigma^\bullet:\K\La 1\Ra^\bullet\otimes_\Z L^\bullet\to M^\bullet
$$
as follows. For every $n\in\Z$, the morphism
$\sigma^n:L^n\oplus L^n\oplus L^{n+1}\to M^n$ restricts to
$\phi^n$ (resp. $\psi^n$, resp. $s^n$) on the first (resp.
second, resp. third) summand. Clearly
\set\begin{equation}\label{eq_from-sigma}
\phi^\bullet=\sigma^\bullet\circ(\iota^\bullet_0\otimes_\Z L^\bullet)
\qquad
\psi^\bullet=\sigma^\bullet\circ(\iota^\bullet_1\otimes_\Z L^\bullet).
\end{equation}
Conversely, the datum of a morphism
$\sigma^\bullet:\K\La 1\Ra^\bullet\otimes_\Z L^\bullet\to M^\bullet$
yields a homotopy from
$\sigma^\bullet\circ(\iota^\bullet_1\otimes_\Z L^\bullet)$
to $\sigma^\bullet\circ(\iota^\bullet_0\otimes_\Z L^\bullet)$.

(iii)\ \
If $\cA$ is abelian, the functor $-\boxtimes_\Z-$ extends
to the category $\sC(\Z\Mod_\mathrm{fg})\times\sC(\cA)$,
and if $\cA$ is also cocomplete, we can even extend it
to the whole of $\sC(\Z\Mod)\times\sC(\cA)$ and we have
as well a corresponding functor $-\otimes_\Z-$ on this
category (see \eqref{subsec_mixed-tensors} and remark
\ref{rem_mixed-tensors}(ii)).

(iv)\ \
Let $M^\bullet$ and $N^\bullet$ be any two complexes of
abelian groups, and $f_0^\bullet,f_1^\bullet:M^\bullet\to N^\bullet$
two homotopy equivalent morphisms. If $K^\bullet$ is any
object of $\sC(\cA)$ and the tensor product
$M^\bullet\otimes_\Z K^\bullet$, $N^\bullet\otimes_\Z K^\bullet$
are both representable in $\sC(\cA)$, then the induced
morphisms $f_0^\bullet\otimes_\Z K^\bullet$ and
$f_1^\bullet\otimes_\Z K^\bullet$ are homotopy equivalent.
Indeed, by (ii), we have a morphism
$\tau^\bullet:\K\La 1\Ra^\bullet\otimes_\Z M^\bullet\to N^\bullet$
such that $\tau^\bullet\circ(\iota^\bullet_i\otimes_\Z\one_M)=
f^\bullet_i$ for $i=0,1$; by example \ref{ex_monoidal}(ii)
the tensor product $\tau^\bullet\otimes_\Z K^\bullet$ is
naturally identified with a morphism
$\tau_K^\bullet:
\K\La 1\Ra^\bullet\otimes_\Z(M^\bullet\otimes_\Z K^\bullet)\to
N^\bullet\otimes_\Z K^\bullet$ such that
$\tau^\bullet_K\circ(\iota^\bullet_i\otimes\one_{M\otimes_\Z K})=
f_i^\bullet\otimes_\Z K^\bullet$ for $i=0,1$, whence the
claim, again by (ii).

(v)\ \
Let $P^\bullet$ be any bounded complex of free abelian
groups of finite rank. Notice the natural identification
$$
(P^\bullet\otimes_\Z K^\bullet)^o\isom
P^{\vee\bullet}\otimes_\Z(K^o)^\bullet
\qquad
\text{in $\sC(\cA^o)$}
$$
deduced from example \ref{ex_dual-is-opposite} and
remark \ref{rem_homotopies}(iii).
\end{remark}

\begin{proposition}\label{prop_univ-hot}
Let $\cA$ be any additive category, and $F:\sC(\cA)\to\cB$
any functor. The following conditions are equivalent :
\begin{enumerate}
\alphaenu
\item
$F$ factors uniquely through the natural functor
$\sC(\cA)\to\Hot(\cA)$.
\item
$F\phi^\bullet$ is an isomorphism, for every homotopy
equivalence $\phi^\bullet$ in $\sC(\cA)$.
\end{enumerate}
\end{proposition}
\begin{proof} Clearly (a)$\Rightarrow$(b), hence, suppose
that (b) holds, consider two chain homotopic morphisms
$\phi^\bullet,\psi^\bullet:L^\bullet\to M^\bullet$ in $\sC(\cA)$,
and define $\sigma^\bullet$ as in remark
\ref{rem_chain-homotopies}(ii); in view of
\eqref{eq_from-sigma}, assertion (a) will follow, once we
know that
$$
F(\iota^\bullet_0\otimes_\Z L^\bullet)=
F(\iota^\bullet_1\otimes_\Z L^\bullet).
$$
To show the latter identity, let
$$
p^\bullet:\K\La 1\Ra^\bullet\to\Z[0]
$$
be the morphism such that $p^0:\Z\oplus\Z\to\Z$ is the
addition map : $(a,b)\mapsto a+b$ for every $a,b\in\Z$.
Clearly
$$
p^\bullet\circ\iota_0^\bullet=\one_{\Z[0]}=
p^\bullet\circ\iota_1^\bullet
$$
so we are easily reduced to showing that
$F(p^\bullet\otimes_\Z L^\bullet)$ is an isomorphism
in $\cB$. By assumption (b), the latter in turn will
follow, once we have shown that
$p^\bullet\otimes_\Z L^\bullet$ is a homotopy equivalence.
Thus, it suffices  to exhibit a morphism
$j^\bullet:\Z[0]\to\K\La 1\Ra^\bullet$ such that
\set\begin{equation}\label{eq_remove-this}
p^\bullet\circ j^\bullet=\one_{\Z[0]}
\end{equation}
and $(j^\bullet\circ p^\bullet)\otimes_\Z L^\bullet$
is homotopy equivalent to $\one_{\K\La 1\Ra\otimes_\Z L}$.
Set $\K\La 1,1\Ra^\bullet:=
\K\La 1\Ra^\bullet\otimes_\Z\K\La 1\Ra^\bullet$; by remark
\ref{rem_chain-homotopies}(ii), the datum of a homotopy from
$(j^\bullet\circ p^\bullet)\otimes_\Z L^\bullet$
to $\one_{\K\La 1\Ra\otimes_\Z L}$ is equivalent to that of
a morphism
$$
\tau^\bullet:
\K\La 1,1\Ra^\bullet\otimes_\Z L^\bullet
\to\K\La 1\Ra^\bullet\otimes_\Z L^\bullet
$$
such that
$$
\tau^\bullet\circ
(\iota_0^\bullet\otimes_\Z\K\La 1\Ra^\bullet\otimes_\Z L^\bullet)=
(j^\bullet\circ p^\bullet)\otimes_\Z L^\bullet
\quad\text{and}\quad
\tau^\bullet\circ
(\iota_1^\bullet\otimes_\Z\K\La 1\Ra^\bullet\otimes_\Z L^\bullet)=
\one_{\K\La 1\Ra\otimes_\Z L}.
$$
Notice that \eqref{eq_remove-this} is fulfilled with
$j^\bullet:=\iota^\bullet_0$; we are then
further reduced to showing the existence of a morphism
$$
t^\bullet:\K\La 1,1\Ra^\bullet\to\K\La 1\Ra^\bullet
$$
such that
$$
t^\bullet\circ
(\iota_0^\bullet\otimes_\Z\K\La 1\Ra^\bullet)=
\iota_0^\bullet\circ p^\bullet
\quad\text{and}\quad
t^\bullet\circ
(\iota_1^\bullet\otimes_\Z\K\La 1\Ra^\bullet)=\one_{\K\La 1\Ra}.
$$
To this aim, let $e_0,e_1$ be the canonical basis
of $\K\La 1\Ra^0$, so that $d^{-1}(1)=e_0-e_1$. Then :
$$
\K\La 1,1\Ra^k
\ \text{is spanned by}\
\left\{\begin{array}{ll}
1\otimes 1 & \text{for $k=-2$} \\
(e\otimes e_i,e_i\otimes e~|~i=0,1) & \text{for $k=-1$} \\
(e_i\otimes e_j~|~i,j=0,1) & \text{for $k=0$}
\end{array}\right.
$$
with differentials given by the rules :
$$
\begin{aligned}
d^{-2}(1\otimes 1):=&\,(e_0-e_1)\otimes 1-1\otimes(e_0-e_1) \\
d^{-1}(1\otimes e_i):=&\,(e_0-e_1)\otimes e_i \\
d^{-1}(e_i\otimes 1):=&\,e_i\otimes(e_0-e_1).
\end{aligned}
$$
We define a morphism $t^\bullet$ as sought, by the rules :
$$
\begin{aligned}
t^{-2}(1\otimes 1):=&\,0 \\
t^{-1}(1\otimes e_0)=t^{-1}(e_0\otimes 1):=&\,0 \\
t^{-1}(e_1\otimes 1)=t^{-1}(1\otimes e_1):=&\,1 \\
t^0(e_i\otimes e_0):=&\,e_0 & & \text{($i=0,1)$} \\
t^0(e_i\otimes e_1):=&\, e_i & & \text{($i=0,1)$}.
\end{aligned}
$$
A direct inspection shows that $t^\bullet$ satisfies
the required identities, and concludes the proof.
\end{proof}

\sset\subsubsection{}\label{subsec_hot-biadditive}
In the situation of \eqref{subsec_biadditive}, let
$A\in\Ob(\cA)$, $A'\in\Ob(\cA')$ and let $P$ be any
free abelian group of finite rank. With the notation
of \eqref{subsec_mixed-tensors}, it is easily seen
that there exist natural isomorphisms in $\cA''$
\set\begin{equation}\label{eq_triple-tens}
B(P\otimes_\Z A,A')\isom P\otimes_\Z B(A,A')\isom
B(A,P\otimes_\Z A').
\end{equation}
Now, let ``$?$'' be either $+$ or $-$, and $K^\bullet$
(resp. $L^\bullet$) any object of $\sC^?(\cA)$ (resp.
of $\sC^?(\cA')$); let also $P^\bullet$ be any bounded
complex of free abelian groups of finite rank. The maps
\eqref{eq_triple-tens} assemble to a natural isomorphism
of triple complexes
$$
\omega^{ijk}:B(P^i\otimes_\Z K^j,L^k)\isom
P^i\otimes_\Z B(K^j,L^k)
\qquad
i,j,k\in\Z
$$
whence, taking into account the discussion of
\eqref{subsec_triple-complex}, an isomorphism of
total complexes
$$
\omega^\bullet_{K,L,P}:
B^\bullet_?(P^\bullet\otimes_\Z K^\bullet,L^\bullet)
\isom P^\bullet\otimes_\Z B^\bullet_?(K^\bullet,L^\bullet)
$$
(notation of \eqref{subsec_cplx-mixed-tens}). Likewise,
we get a natural isomorphism
$$
\sigma^\bullet_{P,K,L}:
B^\bullet_?(K^\bullet,P^\bullet\otimes_\Z L^\bullet)
\isom P^\bullet\otimes_\Z B^\bullet_?(K^\bullet,L^\bullet).
$$
(Notice that $\sigma^\bullet_{P,K,L}$ involves
composition with the flip operator
$\sC(\mathsf{fl}_\cA):\sC_3(\cA)\isom\sC_3(\cA)$,
hence, in each degree $p\in\Z$, it is the direct
sum of morphisms $(-1)^{ij}\sigma^{ijk}$, for
all $i,j,k\in\Z$ such that $i+j+k=p$, and where
$\sigma^{ijk}$ is a natural isomorphism as in
\eqref{eq_triple-tens}).

Especially, take $P^\bullet:=\K\La 1\Ra^\bullet$; it
is easily seen that $\omega^\bullet_{K,L,\K\La 1\Ra}$
identifies the morphism 
$$
B^\bullet_?(\iota_i^\bullet\otimes_\Z\one_K,L^\bullet):
B^\bullet_?(K^\bullet,L^\bullet)\to
B^\bullet_?(\K\La 1\Ra^\bullet\otimes_\Z K^\bullet,L^\bullet)
$$
with the morphism
$$
\iota_i^\bullet\otimes_\Z B^\bullet_?(K^\bullet,L^\bullet):
B^\bullet_?(K^\bullet,L^\bullet)\to
\K\La 1\Ra^\bullet\otimes_\Z B^\bullet_?(K^\bullet,L^\bullet).
$$
Hence, suppose that
$\phi^\bullet_0,\phi_1^\bullet:L^\bullet\to M^\bullet$
are two homotopy equivalent morphisms in $\sC^?(\cA')$;
taking into account remark \ref{rem_chain-homotopies}(ii),
we deduce that the induced morphisms
$B^\bullet_?(K^\bullet,\phi_i^\bullet)$ (for $i=0,1$) are
also homotopy equivalent. A similar argument shows that
the functor
$$
\sC^?(\cA)\to\sC^?(\cA'')
\quad :\quad
K^\bullet\mapsto B^\bullet_?(K^\bullet,L^\bullet)
$$
likewise preserves homotopy equivalences. We conclude
that the functor $B^\bullet_?$ descends to the homotopy
category, to give a functor
$$
B^\bullet_?:\Hot^?(\cA)\times\Hot^?(\cA')\to\Hot^?(\cA'')
$$
which, in case $\cA''$ is cocomplete (resp. complete), can
even be extended to a functor
$$
B^\bullet_\oplus:\Hot(\cA)\times\Hot(\cA')\to\Hot(\cA'')
\qquad
\text{(resp.\ $B^\bullet_\Pi:\Hot(\cA)\times\Hot(\cA')\to\Hot(\cA'')$)}.
$$

\begin{example}\label{ex_hot-morphisms}
(i)\ \
Consider any abelian tensor category
$(\cA,\otimes,\Phi,\Psi)$. Following
\eqref{subsec_hot-biadditive}, the tensor product
of example \ref{ex_monoidal}(i) descends to a
biadditive functor
$$
\otimes:\Hot^?(\cA)\times\Hot^?(\cA)\to\Hot^?(\cA)
$$
for ``$?$'' equal to either $+$ or $-$. If $\cA$
is cocomplete, this functor is even defined on the
whole of $\Hot(\cA)$.

(ii)\ \
The constructions of \eqref{subsec_hot-biadditive}
can also be applied to the functor $\Hom^\bullet_\cA$
of example \ref{ex_Hom-complex-hots}(i) : taking into
account remark \ref{rem_homotopies}(iii), we deduce
that the functor $\Hom^\bullet_\cA$ descends to a
biadditive functor :
$$
\rHot_\cA^\bullet:\Hot(\cA)\times\Hot(\cA)^o\to\Hot(\Z\Mod)
\quad :\quad
(L^\bullet,K^\bullet)\mapsto\Hom_\cA^\bullet(K^\bullet,L^\bullet).
$$
\end{example}

\subsection{Filtered complexes and spectral sequences}
\label{sec_filter-cpx}
Let $\cA$ be any abelian category; a {\em filtered object}
of $\cA$ is a datum
$$
(A,\Fil^\bullet A)
$$
consisting of an object $A$ of $\cA$ and a system of
subobjects $(\Fil^pA~|~p\in\Z)$ of $A$, with
$$
\Fil^{p+1}A\subset\Fil^pA
\qquad
\text{for every $p\in\Z$}.
$$
If $(B,\Fil^\bullet B)$ is any other filtered object of
$\cA$, a morphism of filtered objects
\set\begin{equation}\label{eq_fil-map}
(A,\Fil^\bullet A)\to
(B,\Fil^\bullet B)
\end{equation}
is just a morphism $u:A\to B$ in $\cA$ which restricts
to a morphism
$$
\Fil^pu:\Fil^pA\to\Fil^pB
\qquad
\text{for every $p\in\Z$}.
$$
With these morphisms, clearly the filtered objects
of $\cA$ form a category
$$
\sFil.\cA.
$$
Moreover, any functor $F:\cA\to\cB$ induces a functor
$$
\sFil.F:\sFil.\cA\to\sFil.\cB
$$
that assigns to any filtered object $(A,\Fil^\bullet A)$
of $\cA$ the filtered object $(FA,\Fil^\bullet FA)$ such
that $\Fil^kFA:=\Img(F(\Fil^kA)\to FA)$ for every $k\in\Z$.

\sset\subsubsection{}\label{subsec_finite-fils}
To any object $(A,\Fil^\bullet A)$ of $\sFil.\cA$ we
attach its {\em associated graded object}, which is
the system of objects of $\cA$ :
$$
(\gr^nA~|~n\in\Z)
\qquad\text{where}\qquad
\gr^nA:=\Fil^nA/\Fil^{n+1}A
\quad
\text{for every $n\in\Z$}.
$$
Notice that any morphism $u(A,\Fil^\bullet)\to(B,\Fil^\bullet)$
of filtered objects of $\cA$ induces a system of morphisms
in $\cA$ :
$$
\gr^pu:\gr^pA\to\gr^pB
\qquad
\text{for every $p\in\Z$}.
$$
We say that the filtration $\Fil^\bullet A$ of a
filtered object $(A,\Fil^\bullet A)$ is {\em bounded above}
(resp. {\em bounded below}) if there exists
$N\in\Z$ such that $\Fil^NA=A$ (resp. such that $\Fil^NA=0$).
We say that the filtration $\Fil^\bullet$ is {\em finite}
if it is bounded above and below. Clearly, if
$\Fil^\bullet A$ is bounded above (resp. below), then
$\gr^nA=0$ whenever $-n$ (resp. $n$) is sufficiently large. 

\begin{definition}\label{def_spectral-seq}
Let $\cA$ be an abelian category, and $r_0\in\N$ any
integer.

\begin{enumerate}
\item
A {\em (cohomological) $r_0$-spectral sequence}
in $\cA$ is a datum
$$
((E^{pq}_r,d_r^{pq},\beta^{pq}_r)~|~p,q,r\in\Z,\ r\geq r_0)
$$
consisting of :
\begin{itemize}
\item
objects $E^{pq}_r$ of $\cA$ and morphisms
$d^{pq}_r:E^{pq}_r\to E^{p+r,q-r+1}_r$ in $\cA$, such that :
\set\begin{equation}\label{eq_complex-cond}
d^{p+r,q-r+1}_r\circ d_r^{pq}=0
\qquad
\text{for every $p,q,r\in\Z$ with $r\geq r_0$}.
\end{equation}
\item
isomorphisms
$$
\beta^{pq}_r:
\Ker\,d_r^{pq}/\Img\,d_r^{p-r,q+r-1}\isom E_{r+1}^{pq}
\qquad
\text{for every $p,q,r\in\Z$ with $r\geq r_0$}.
$$
\end{itemize}
\item
Let $(E^{\bullet\bullet}_\bullet,d^{\bullet\bullet}_{E,\bullet},
\beta^{pq}_{E,\bullet})$ and
$(F^{\bullet\bullet}_\bullet,,d^{\bullet\bullet}_{F,\bullet},
\beta^{pq}_{F,\bullet})$ be two $r_0$-spectral sequences
in $\cA$. A {\em morphism of $r_0$-spectral sequences}
$$
u^{\bullet\bullet}_\bullet:E^{\bullet\bullet}_\bullet\to
F^{\bullet\bullet}_\bullet
$$
is a system $(u^{pq}_r:E^{pq}_r\to F^{pq}_r~|~p,q,r\in\Z,\ r\geq r_0)$
of morphisms of $\cA$, such that the diagrams
$$
\xymatrix{
E^{pq}_r \ar[r]^-{d^{pq}_{E,r}} \ar[d]_{u^{pq}_r} &
E^{p+r,q-r+1} \ar[d]^{u^{p+r,q-r+1}_r} &
E^{pq}_r \ar[d]^{u^{pq}_r} &
\Ker\,d^{pq}_{E,r} \ar[r] \ar[l] &
\Ker\,d^{pq}_{E,r}/\Img\,d^{p-r,q+r-1}_{E,r}
\ar[r]^-{\beta_{E,r}^{pq}} & E^{pq}_{r+1} \ar[d]_{u^{pq}_{r+1}} \\
F^{pq}_r \ar[r]^-{d^{pq}_{F,r}} & F^{p+r,q-r+1} &
F^{pq}_r &
\Ker\,d^{pq}_{F,r} \ar[r] \ar[l] &
\Ker\,d^{pq}_{E,r}/\Img\,d^{p-r,q+r-1}_{E,r}
\ar[r]^-{\beta_{F,r}^{pq}} & F^{pq}_{r+1}
}$$
commute for every $p,q,r\in\Z$ with $r\geq r_0$
(where the unmarked arrows are the natural monomorphisms
and epimorphisms). Obviously, the $r_0$-spectral sequences
of $\cA$ and their morphisms form a category, which
we denote by
$$
\SpSeq_{r_0}(\cA)
$$
\end{enumerate}
\end{definition}

\sset\subsubsection{}\label{subsec_decalage-spseq}
Let $r\in\N$ be any integer. Clearly, for every integer
$s\geq r$ there is a natural functor
$$
\SpSeq_r(\cA)\to\SpSeq_s(\cA)
$$
that simply forgets the terms $E^{pq}_k$ with $r\leq k<s$.
One may also define a natural functor form $s$-spectral
sequences to $r$-spectral sequences. Namely, let
$E^{\bullet\bullet}_\bullet$ be any $(r+1)$-spectral
sequence of $\cA$. The {\em d\'ecalage} of
$E^{\bullet\bullet}$ is the $r$-spectral sequence
$(\Dec(E)^{\bullet\bullet},d^{\bullet\bullet}_{\Dec(E)})$
given by the rule
$$
\Dec(E)^{p,q-p}_s:=E^{p+q,-p}_{s+1}
\qquad
d_{\Dec(E),s}^{p,q-p}:=d^{p+q,-p}_{E,s+1}
\qquad
\text{for every $p,q,s\in\Z$ with $s\geq r$}
$$
and whose isomorphisms $\beta_{\Dec(E)}^{\bullet\bullet}$
are deduced from the corresponding isomorphisms for
$E^{\bullet\bullet}$, in the obvious way. Clearly, this
rule extends to a natural functor
$$
\Dec:\SpSeq_{r+1}(\cA)\to\SpSeq_r(\cA).
$$

\sset\subsubsection{}\label{subsec_spectral-abut}
Let $E^{\bullet\bullet}_\bullet$ be any $r_0$-spectral sequence;
with the notation of definition \ref{def_spectral-seq}
we define, by induction on $k-r$, two systems of subobjects
$$
(Z_k(E^{pq}_r),B_k(E^{pq}_r)\subset E^{pq}_r~|~
p,q,r,k\in\Z,\ k>r\geq r_0)
$$
as follows. First, we set
$$
Z_{r+1}(E^{pq}_r):=\Ker\,d^{pq}_r
\qquad
B_{r+1}(E^{pq}_r):=\Img\,d^{p-r,q+r-1}_r
\qquad
\text{for every $p,q,r\in\Z$ with $r\geq r_0$}.
$$
Next, say that $i\geq 2$ is any integer, and suppose
that both $Z_k(E^{pq}_r)$ and $B_k(E^{pq}_r)$ have
already been defined, for every $p,q,r,k\in\Z$ with
$k>r\geq r_0$ and $k-r<i$; we let
$$
Z_{r+i}(E^{pq}_r):=\pi_{pqr}^{-1}Z_{r+i}(E^{pq}_{r+1})
\qquad
B_{r+i}(E^{pq}_r):=\pi^{-1}_{pqr}B_{r+i}(E^{pq}_{r+1})
$$
where $\pi_{pqr}$ is the composition
$$
Z_{r+1}(E^{pq}_r)\to Z_{r+1}(E^{pq}_r)/B_{r+1}(E^{pq}_r)
\xrightarrow{\ \beta^{pq}_r\ }E^{pq}_{r+1}. 
$$
Then it is clear that the system
$(\beta^{pq}_r~|~p,q,r\in\Z,\ r\geq r_0)$ induces
isomorphisms
$$
\beta^{pq}_{rk}:Z_k(E^{pq}_r)/B_k(E^{pq}_r)\isom E^{pq}_k
\qquad
\text{for every $p,q,r,k\in\Z$ with $k>r\geq r_0$}
$$
and for every $p,q,r,k,n\in\Z$ with $k,n>r\geq r_0$
we have inclusions
$$
B_k(E^{pq}_r)\subset B_{k+1}(E^{pq}_r)\subset
Z_{n+1}(E^{pq}_r)\subset Z_n(E^{pq}_r)\subset E^{pq}_r.
$$

\begin{definition}\label{def_abutment}
In the situation of \eqref{subsec_spectral-abut}, we define :

(i)\ \
an {\em abutment} for the spectral sequence
$E^{\bullet\bullet}_\bullet$ to be a datum consisting of
\begin{itemize}
\item
a system of subobjects
$(B_\infty(E^{pq}_{r_0}),Z_\infty(E_{r_0}^{pq})
\subset E^{pq}_{r_0}~|~p,q\in\Z)$ such that
$$
B_k(E^{pq}_{r_0})\subset B_\infty(E_{r_0}^{pq})
\subset Z_\infty(E_{r_0}^{pq})\subset Z_k(E^{pq}_{r_0})
\qquad
\text{for every $p,q,k\in\Z$ with $k>r_0$}
$$
\item
a system of filtered objects
$((E^n_\infty,\Fil^\bullet E^n_\infty)~|~n\in\Z)$  of $\cA$
\item
a system of isomorphisms in $\cA$
$$
\beta_\infty^{pq}:E_\infty^{pq}:=
Z_\infty(E_{r_0}^{pq})/B_\infty(E^{pq}_{r_0})
\isom\gr^pE^{p+q}_\infty
\qquad
\text{for every $p,q\in\Z$}. 
$$
\end{itemize}
We summarize these conditions via the traditional
notation :
$$
E^{pq}_r\Rightarrow E^{p+q}_\infty.
$$
For any $n\in\Z$, we say that $E^{\bullet\bullet}_\bullet$ is
{\em convergent in degree $n$}, if for every $p,q\in\Z$
with $p+q=n$ there exists $k\in\N$ such that
$B_k(E^{pq}_{r_0})=B_\infty(E^{pq}_{r_0})$ and
$Z_k(E^{pq}_{r_0})=Z_\infty(E^{pq}_{r_0})$. We say that
$E^{\bullet\bullet}_\bullet$ is {\em convergent}, if it is
convergent in all degrees.

(iii)\ \
Let  $(E^{\bullet\bullet}_\bullet,B_\infty(E_{r_0}^{\bullet\bullet}),
Z_\infty(E_{r_0}^{\bullet\bullet}),E_\infty^\bullet,\beta^{pq}_{E,\infty})$
and $(F^{\bullet\bullet}_\bullet,B_\infty(F_{r_0}^{\bullet\bullet}),
Z_\infty(F_{r_0}^{\bullet\bullet}),F_\infty^\bullet,\beta^{pq}_{F,\infty})$
be two $r_0$-spectral sequences in $\cA$ with abutments.
A {\em morphism of $r_0$-spectral sequences with abutments}
is a pair
$$
u^{\bullet\bullet}_\bullet:E^{\bullet\bullet}_\bullet\to
F^{\bullet\bullet}_\bullet
\qquad
(v^n:(E^n_\infty,\Fil^\bullet E^n_\infty)\to
(F^n_\infty,\Fil^\bullet F^n_\infty)~|~n\in\Z)
$$
such that :
\begin{itemize}
\item
$u^{\bullet\bullet}_\bullet$ is a morphism of
$r_0$-spectral sequences.
\item
$v^n$ is a morphism of filtered objects of $\cA$,
for every $n\in\Z$.
\item
$u^{pq}_{r_0}(B_\infty(E^{pq}_{r_0}))\subset
B_\infty(F_{r_0}^{pq})$ and
$u^{pq}_{r_0}(Z_\infty(E^{pq}_{r_0}))\subset
Z_\infty(F_{r_0}^{pq})$
for every $p,q\in\Z$.
\item
The following diagrams commute, for every $p,q\in\Z$ :
$$
\xymatrix{ E^{pq}_{r_0} \ar[d]_{u^{pq}_{r_0}} &
\ar[l] Z_\infty(E^{pq}_{r_0}) \ar[r] &
E_\infty^{pq} \ar[rr]^-{\beta_{E,\infty}^{pq}} & &
\gr^pE^{p+q}_\infty \ar[d]^{\gr^pv^{p+q}} \\
F^{pq}_{r_0} & \ar[l] Z_\infty(F^{pq}_{r_0}) \ar[r] &
E_\infty^{pq} \ar[rr]^-{\beta_{F,\infty}^{pq}} & &
\gr^pF^{p+q}_\infty
}$$
(where the unmarked arrows are the natural epimorphisms
and monomorphisms).
\end{itemize}
\end{definition}

\sset\subsubsection{}\label{subsec_simply-forget}
For every $r\in\N$, we denote by
$$
\SpSeq_r^\infty(\cA)
$$
the category of $r$-spectral sequences with abutments (and
their morphisms, as in definition \ref{def_abutment}(ii)).
For every integer $s>r$ there is an obvious functor
$$
\SpSeq^\infty_r(\cA)\to\SpSeq_s^\infty(\cA)
$$
that simply forgets the terms $E^{pq}_k$ with $r\leq k<s$,
and replaces the terms $B_\infty(E^{pq}_r),Z_\infty(E_r^{pq})$
with their images
$\beta^{pq}_{rs}(B_\infty(E^{pq}_r)),\beta^{pq}_{rs}(Z_\infty(E_r^{pq})$
in $E^{pq}_s$.

\begin{remark}\label{rem_vanish-reproduces}
Let  $r_0\in\N$, $n\in\Z$ and $E^{\bullet\bullet}_\bullet$ any
$r_0$-spectral sequence of $\cA$.

(i)\ \
Notice that if $E^{pq}_r=0$ for some $p,q\in\Z$
and $r\geq r_0$, then $E^{pq}_s=0$ for every $s\geq r$. 

(ii)\ \
We say $E^{\bullet\bullet}_\bullet$ is {\em bounded above}
(resp. {\em bounded below}) in degree $n$, if there exist
$r,N\in\N$ with $r\geq r_0$, such that
$$
E_r^{pq}=0
\qquad
\text{for all $p,q\in\Z$ such that $p+q=n$ and $q>N$
(resp. $p>N$)}.
$$
If $E^{\bullet\bullet}_\bullet$ is bounded above in degree $n+1$
and bounded below in degree $n-1$, notice that for every
$p,q\in\Z$ with $p+q=n$ there exists an integer $r\geq r_0$
such that both $d^{pq}_s$ and $d_s^{p-s,q+s-1}$ are zero morphisms
for every $s\geq r$, in which case $\beta^{pq}_s$ is an isomorphism
$E_{s+1}^{pq}\isom E^{pq}_s$, for every  $s\geq r$, and for every
$p,q,s,k\in\Z$ with $k>s\geq r$ and $p+q=n$ we have :
$$
Z_k(E^{pq}_s)=Z_{s+1}(E^{pq}_s)
\qquad
B_k(E^{pq}_s)=B_{s+1}(E^{pq}_s).
$$
\end{remark}

\begin{definition}
Let $\cA$ be any abelian category, and $n\in\Z$.

(i)\ \
A {\em filtered complex} $(K^\bullet,\Fil^\bullet K^\bullet)$
of $\cA$ is just an object of the category $\sFil.\sC(\cA)$.

(ii)\ \
We say that $\Fil^\bullet K^\bullet$ is {\em bounded above}
(resp. {\em bounded below}, resp. {\em finite}) in degree $n$,
if the filtration $\Fil^\bullet K^n$ on the object $K^n$
of $\cA$ is bounded above (resp. bounded below, resp. finite).
We say that $\Fil^\bullet K^\bullet$ is {\em bounded}, if it is
finite in all degrees (see \eqref{subsec_finite-fils}).
\end{definition}

\begin{remark}\label{rem_fil-tot}
Likewise, we may define a {\em filtered double complex}
as an object of the category $\sFil.\sC_2(\cA)$. Then,
the flip and diagonal functors of \eqref{subsec_flip-diag}
can be upgraded to functors
$$
\sFil.\sC_2(\cA)\to\sFil.\sC_2(\cA)
\qquad
\sFil.\sC_2(\cA)\to\sFil.\sC(\cA).
$$
If all coproducts (resp. all products) are representable in
$\cA$, the total complex functor $\Tot^\oplus$ (resp. $\Tot^\Pi$)
induces a corresponding functor
$$
\sFil.\Tot^\oplus:\sFil.\sC_2(\cA)\to\sFil.\sC(\cA)
\qquad
\text{(resp.\  $\sFil.\Tot^\Pi:\sFil.\sC_2(\cA)\to\sFil.\sC(\cA)$}
$$
(see \eqref{sec_filter-cpx}) which we shall often denote
just by $\Tot$ (resp. by $\Tot^\Pi$).
\end{remark}

\sset\subsubsection{}\label{subsec_fil-spec-seq}
Any filtered complex of $\cA$ determines a spectral
sequence of $\cA$, whose terms are defined as follows.
For every $p,q\in\Z$ and $r\in\N$, set
$$
\begin{aligned}
Z(K)_r^{pq}:=\,& \Ker\,
(\bar d:(\Fil^pK^\bullet)^{p+q}\to
(\Fil^pK^\bullet)^{p+q+1}/(\Fil^{p+r}K^\bullet)^{p+q+1}) \\
D(K)^{pq}_r:=\,&
(\Fil^{p+1}K^\bullet)^{p+q}+d_K^{p+q-1}((\Fil^{p-r+1}K^\bullet)^{p+q-1})
\subset(\Fil^{p-r+1}K^\bullet)^{p+q} \\
B(K)^{pq}_r:=\,& D(K)_r^{pq}\cap Z(K)^{pq}_r \\
E(K)^{pq}_r:=\,& Z(K)^{pq}_r/B(K)^{pq}_r
\end{aligned}
$$
where $\bar d$ is the morphism in $\cA$ induced
by the differential
$$
d^{p+q}_{\Fil^pK}:
(\Fil^pK^\bullet)^{p+q}\to(\Fil^pK^\bullet)^{p+q+1}
$$
of the complex $\Fil^pK^\bullet$. Notice that
$$
(\Fil^{p+r}K^\bullet)^{p+q}\subset Z(K)_r^{pq}
\subset(\Fil^pK^\bullet)^{p+q}
\qquad
\text{for every $p,q\in\Z$ and $r\in\N$}.
$$

\begin{lemma}\label{lem_notice-that}
With the notation of \eqref{subsec_fil-spec-seq}, we have :
$$
B(K)_r^{pq}=Z(K)_r^{p+1,q-1}+
d_K^{p+q-1}(Z(K)_{r-1}^{p-r+1,q+r-2})
\qquad
\text{for every $p,q\in\Z$ and $r\geq 1$}.
$$
\end{lemma}
\begin{proof} Notice first that
$$
D(K)_r^{pq}\cap(\Fil^pK^\bullet)^{p+q}=(\Fil^{p+1}K^\bullet)^{p+q}
+(\Img\,d_{\Fil^{p-r+1}K}^{p+q-1}\cap(\Fil^pK^\bullet)^{p+q})
$$
and
$$
\Img\,d_{\Fil^{p-r+1}K}^{p+q-1}\cap(\Fil^pK^\bullet)^{p+q}=
d^{p+q-1}_K(Z(K)_{r-1}^{p-r+1,q+r-2}).
$$
On the other hand, we have :
$$
d^{p+q}_K(D(K)^{pq}_r)\cap(\Fil^{p+r}K^\bullet)^{p+q+1}=
\Img\,d_{\Fil^{p+1}K}^{p+q}\cap(\Fil^{p+r}K^\bullet)^{p+q+1}=
d^{p+q}_K(Z(K)_r^{p+1,q-1}).
$$
The lemma follows easily, by comparing these identities.
\end{proof}

\sset\subsubsection{}\label{subsec_spec-seq-fil-cpx}
It is easily seen that $d^{p+q}_{\Fil^pK}$ restricts
to a morphism
\set\begin{equation}\label{eq_beg-to-differ}
D(K)^{pq}_r+Z(K)_r^{pq}\to Z(K)_r^{p+r,q-r+1}
\end{equation}
which sends $D(K)^{pq}_r$ into $B(K)_r^{p+r,q-r+1}$, and
therefore induces a morphism
$$
d(K)^{pq}_r:E(K)^{pq}_r\to E(K)^{p+r,q-r+1}_r
\qquad
\text{for every $p,q\in\Z$ and $r\in\N$}.
$$
Clearly $d(K)^{p+r,q-r+1}_r\circ d(K)_r^{pq}=0$.
Furthermore, a simple inspection shows that
$$
\Ker\,d(K)_r^{pq}=Z(K)^{pq}_{r+1}/(Z(K)^{pq}_{r+1}\cap D(K)^{pq}_r)
\qquad
\Img\,d(K)_r^{p-r,q+r-1}=D(K)^{pq}_{r+1}/D(K)^{pq}_r
$$
whence a natural isomorphism
$$
E(K)_{r+1}^{pq}\isom\Ker\,d(K)_r^{pq}/\Img\,d(K)_r^{p-r,q+r-1}
\qquad
\text{for every $p,q\in\Z$ and $r\in\N$}
$$
so the system
$(E(K)^{\bullet\bullet}_\bullet,d(K)^{\bullet\bullet}_\bullet)$
is indeed a spectral sequence. 

\begin{remark}\label{rem_low-terms}
(i)\ \
A simple inspection shows that
$$
E(K)_0^{pq}=\gr^pK^{p+q}
\qquad\text{and}\qquad
d(K)_0^{pq}=d_{\gr^pK}^{p+q}
\qquad
\text{for every $p,q\in\Z$}
$$
so the resulting complex $(E(K)_0^{p,\bullet},d(K)^{p,\bullet})$
is just $\gr^pK^\bullet[q]$, for every $p,q\in\Z$.

(ii)\ \
In view of (i), we have a natural isomorphism
$$
E(K)^{pq}_1\isom H^{p+q}(\gr^pK^\bullet)
\qquad
\text{for every $p,q\in\Z$}
$$
which identifies $d(K)_1^{pq}$ with a morphism
$$
H^{p+q}(\gr^pK^\bullet)\to H^{p+q+1}(\gr^{p+1}K^\bullet).
$$
A direct inspection shows that the latter is the
boundary morphism in degree $p+q$ arising (by the
snake lemma) from the short exact of complexes
$$
0\to\gr^{p+1}K^\bullet\to
\Fil^pK^\bullet/\Fil^{p+2}K^\bullet\to\gr^pK^\bullet\to 0.
$$
\end{remark}

\sset\subsubsection{}\label{subsec_Fil-toSpSeq}
With the notation of \eqref{subsec_spectral-abut} and
\eqref{subsec_fil-spec-seq}, it is easily seen that
$$
\begin{aligned}
Z_k(E(K)^{pq}_r)=\, & \Img\,(Z(K)^{pq}_k\to E(K)^{pq}_r)  \\
B_k(E(K)^{pq}_r)=\, & \Img\,(B(K)^{pq}_k\to E(K)^{pq}_r)
\end{aligned}
$$
for every $p,q,r,k\in\Z$ with $k>r\geq 0$. We may then
define a natural abutment for $E(K)^{\bullet\bullet}_\bullet$,
as follows. For every $p,q\in\Z$ we set
$$
\begin{aligned}
Z_\infty(E(K)_0^{pq}):=\, &
\Img(\Ker\,d^{p+q}_{\Fil^pK}\to E(K)^{pq}_0) \\
B_\infty(E(K)_0^{pq}):=\, &
\Img((\Fil^pK^{p+q}\cap\Img\,d^{p+q-1}_K)\to E(K)^{pq}_0)
\end{aligned}
$$
as well as
$$
E(K)_\infty^n:=H^nK^\bullet
\qquad
\Fil^pE(K)^n_\infty:=\Img(H^n(\Fil^pK^\bullet)\to H^n(K^\bullet))
\qquad
\text{for every $p,n\in\Z$}.
$$
By inspecting the definitions we find natural isomorphisms
$$
\beta(K)_\infty^{pq}:E(K)_\infty^{pq}:=
Z_\infty(E(K)_0^{pq})/B_\infty(E(K)^{pq}_0)
\isom\gr^pE(K)^{p+q}_\infty
$$
and the datum
$A(K)^{\bullet\bullet}_\infty:=
(B_\infty(E(K)_0^{pq}),Z_\infty(E(K)_0^{pq}),\beta(K)_\infty^{pq})$
is the sought abutment. Summing up, we get a well defined functor
$$
\sFil.\sC(\cA)\to\SpSeq_0^\infty(\cA)
\qquad
(K^\bullet,\Fil^\bullet K^\bullet)\mapsto
(E(K)^{\bullet\bullet}_\bullet,A(K)^{\bullet\bullet}_\infty).
$$
Clearly, if $\Fil^\bullet K^\bullet$ is bounded in degree $n$,
then $\Fil^\bullet E(K)^n_\infty$ is a finite filtration.

\begin{proposition}\label{prop_convergence-simple}
{\em (i)}\ \
Let $(K^\bullet,\Fil^\bullet K^\bullet)$ be a filtered complex
in $\cA$, and $N\in\N$ such that :
$$
H^{n-1}(K^\bullet/\Fil^{-p}K^\bullet)=0
\qquad\text{and}\qquad
H^{n+1}(\Fil^pK^\bullet)=0
\qquad
\text{for every $p\geq N$}.
$$
Then the $1$-spectral sequence with abutment
$(E(K)^{\bullet\bullet}_\bullet,A(K)^{\bullet\bullet}_\infty)$
is convergent in degree $n$.

{\em(ii)}\ \
Especially, the conclusion of {\em (i)} holds if
$\Fil^\bullet K^\bullet$ is bounded above in degree
$n-1$ and bounded below in degree $n+1$.
\end{proposition}
\begin{proof} Obviously (i)$\Rightarrow$(ii). Hence, let us
assume that the condition of (i) holds, and consider, for every
$r\in\N$ and every $p,q\in\Z$ with $p+q=n$, the composition
$$
H^n(\gr^pK^\bullet)\xrightarrow{\delta}H^{n+1}(\Fil^{p+1}K^\bullet)
\xrightarrow{H^{n+1}(\pi_r)}H^{n+1}(\Fil^{p+1}K^\bullet/\Fil^{p+r}K^\bullet)
$$
where $\delta$ is the boundary morphism attached by
the snake lemma to the short exact sequence of complexes
$0\to\Fil^{p+1}K^\bullet\!\to\!\Fil^pK^\bullet\!\to\!\gr^pK^\bullet\to 0$,
and $\pi_r:\Fil^{p+1}K^\bullet\to\Fil^{p+1}K^\bullet/\Fil^{p+r}K^\bullet$
is the projection. By direct inspection, we see that :
$$
\Ker\,(H^{n+1}(\pi_r)\circ\delta)=Z_r(E_1^{pq})
\qquad\text{and}\qquad
\Ker\,\delta=Z_\infty(E_1^{pq}).
$$
On the other hand, the long exact cohomology sequence
attached to the short exact sequence of complexes
$0\to\Fil^{p+r}K^\bullet\to\Fil^{p+1}K^\bullet\to
\Fil^{p+1}K^\bullet/\Fil^{p+r}K^\bullet\to 0$ shows that
$\Ker\,H^{n+1}(\pi_r)$ is a quotient of
$H^{n+1}(\Fil^{p+r}K^\bullet)$. Thus, $H^{n+1}(\pi_r)$ is a
monomorphism for every $r\in\N$ such that $p+r\geq N$,
and in that case we get $Z_r(E_1^{pq})=Z_\infty(E_1^{pq})$.
Likewise, consider the composition :
$$
H^{n-1}(\Fil^{p-r+1}K^\bullet/\Fil^pK^\bullet)
\xrightarrow{H^{n-1}(j_r)}H^{n-1}(K^\bullet/\Fil^pK^\bullet)
\xrightarrow{\delta'}H^n(\gr^pK^\bullet)
$$
where $\delta'$ is the boundary morphism attached by the
snake lemma to the short exact sequence
$0\to\Fil^pK^\bullet\to K^\bullet\to K^\bullet/\Fil^pK^\bullet\to 0$,
and $j_r:\Fil^{p-r+1}K^\bullet/\Fil^pK^\bullet\to
K^\bullet/\Fil^pK^\bullet$ is the natural inclusion. By direct
inspection we see that
$$
\Img(\delta'\circ H^{n-1}(j_r))=B_r(E_1^{pq})
\qquad\text{and}\qquad
\Img\,\delta'=B_\infty(E^{pq}_1).
$$
On the other hand, the long exact cohomology sequence
attached to the short exact sequence
$0\to\Fil^{p-r+1}K^\bullet/\Fil^pK^\bullet\to
K^\bullet/\Fil^pK^\bullet\to K^\bullet/\Fil^{p-r+1}K^\bullet\to 0$
shows that $\Coker\,H^{n-1}(j_r)$ is a subobject of
$H^{n-1}(K^\bullet/\Fil^{p-r+1}K^\bullet)$. Especially,
$H^{n-1}(j_r)$ is an epimorphism whenever $p-r+1\leq -N$,
and in that case we get $B_r(E_1^{pq})=B_\infty(E_1^{pq})$.
\end{proof}

\sset\subsubsection{}
The constructions of examples \ref{ex_Hom-complex-hots}(i),
\ref{ex_monoidal}(i) admit filtered counterparts.
Namely, suppose that $(\cA,\otimes,\Phi,\Psi)$ is a
tensor abelian category, and let $(A,\Fil^\bullet A)$
and $(B,\Fil^\bullet B)$ be any two objects of $\sFil.\cA$.
Suppose that :
\begin{itemize}
\item
at least one of the filtrations $\Fil^\bullet A$ and
$\Fil^\bullet B$ is bounded
\item
or else, $\cA$ is cocomplete.
\end{itemize}
In this situation we define a filtration
$\Fil^\bullet(A\otimes B)$ on $A\otimes B$, by ruling that
$$
\Fil^p(A\otimes B):=
\sum_{k\in\Z}\Img\,(\Fil^kA\otimes\Fil^{p-k}B\to A\otimes B)
\qquad
\text{for every $p\in\Z$}.
$$
In case $\cA$ is cocomplete, clearly, this rule yields a
functor
$$
\sFil.\cA\times\sFil.\cA\to\sFil.\cA.
$$
Next, suppose that $\cA$ is cocomplete, and let
$(K^\bullet,\Fil^\bullet K^\bullet)$ and
$(L^\bullet,\Fil^\bullet L^\bullet)$ be any two filtered
complexes of $\cA$; we get an induced filtration
$\Fil^\bullet K^\bullet\boxtimes L^\bullet$ on the double
complex $K^\bullet\boxtimes L^\bullet$, and applying
the functor $\Tot$ of remark \ref{rem_fil-tot} we obtain
a filtration on $K^\bullet\otimes L^\bullet$, whence a
functor
$$
\sFil.\sC(\cA)\times\sFil.\sC(\cA)\to\sFil.\sC(\cA)
\qquad
(\Fil^\bullet K,\Fil^\bullet L)\mapsto
\Fil^\bullet K^\bullet\otimes L^\bullet.
$$
Notice that if $\cA$ is not cocomplete, the filtered double
complex $\Fil^\bullet K^\bullet\boxtimes L^\bullet$ is still
well defined, provided at least one of the filtrations
$\Fil^\bullet K^\bullet$ and $\Fil^\bullet L^\bullet$ is
bounded. If moreover both $K^\bullet$ and $L^\bullet$
lie in $\sC^-(\cA)$ or $\sC^+(\cA)$, then also
$\Fil^\bullet K^\bullet\otimes L^\bullet$ is well defined.

\sset\subsubsection{}
Let $(A,\Fil^\bullet A)$ and $(B,\Fil^\bullet B)$ be
any two objects of $\sFil.\cA$. We define a filtration
$\Fil^\bullet\Hom_\cA(A,B)$ on the abelian group
$\Hom_\cA(A,B)$, by ruling that
$$
\Fil^p\Hom_\cA(A,B):=
\{f:A\to B~|~f(\Fil^kA)\subset\Fil^{p+k}B\ 
\text{for every $k\in\Z$}\}
$$
for every $p\in\Z$. Clearly this rule yields a functor
\set\begin{equation}\label{eq_apply-this-fun}
\sFil.\cA\times(\sFil.\cA)^o\to\sFil.\Z\Mod.
\end{equation}
Let now $(K^\bullet,\Fil^\bullet K^\bullet)$ and
$(L^\bullet,\Fil^\bullet L^\bullet)$ be any two filtered
complexes of $\cA$; we define a filtration
$\Fil^\bullet\Hom^{\bullet\bullet}_\cA(K^\bullet,L^\bullet)$
on the double complex of abelian groups
$\Hom^{\bullet\bullet}_\cA(K^\bullet,L^\bullet)$ (see example
\ref{ex_Hom-complex-hots}(i)), by ruling that
$$
\Fil^k\Hom^{p,q}_\cA(K^\bullet,L^\bullet):=
\Fil^k\Hom_\cA(K^{-p},L^q)
\qquad
\text{for every $p,q,k\in\Z$}
$$
(where the right-hand side is the filtered abelian group
obtained by applying the functor \eqref{eq_apply-this-fun}
to the filtered objects $(K^{-p},\Fil^\bullet K^{-p})$
and $(L^q,\Fil^\bullet L^q)$ of $\cA$). After applying
the functor $\Tot^\Pi$ of remark \ref{rem_fil-tot}, we
deduce a filtration on $\Hom_\cA^\bullet(K^\bullet,L^\bullet)$,
whence a functor
$$
\sFil.\sC(\cA)\times\sFil.\sC(\cA)^o\to\sFil.\sC(\Z\Mod)
\qquad
(\Fil^\bullet L^\bullet,\Fil^\bullet K^\bullet)
\mapsto\Fil^\bullet\Hom^\bullet_\cA(K^\bullet,L^\bullet).
$$
Taking into account remark \ref{ex_Hom-complex-hots}(i),
we also see that
$$
\Hom_{\sFil.\sC(\cA)}(K^\bullet,L^\bullet)=
\Fil^0\Hom^\bullet(K^\bullet,L^\bullet)\cap\ker\,d^0_{K,L}.
$$

\begin{definition}\label{def_fil-homotopy}
Let $\cA$ be any abelian category, $k\in\Z$ any integer,
$\phi^\bullet,\psi^\bullet:(K^\bullet,\Fil^\bullet K^\bullet)\to
(L^\bullet,\Fil^\bullet L^\bullet)$ any two morphisms of
filtered complexes of $\cA$.
\begin{enumerate}
\item
A {\em homotopy of order $k$} from $\phi^\bullet$ to
$\psi^\bullet$ is an element
$s^\bullet\in\Fil^k\Hom^{-1}_\cA(K^\bullet,L^\bullet)$
such that $d^{-1}_{K,L}(s^\bullet)=\psi^\bullet-\phi^\bullet$
(notation of example \ref{ex_Hom-complex-hots}(i)).
\item
We write $\phi^\bullet\sim_k\psi^\bullet$ if there
exists a homotopy of order $k$ from $\psi^\bullet$
to $\phi^\bullet$. It is easily seen that if this
is the case, and $\alpha^\bullet:K'{}^\bullet\to K^\bullet$,
$\beta^\bullet:L^\bullet\to L'{}^\bullet$ are any two morphisms,
then $\phi^\bullet\circ\alpha^\bullet\sim_k\psi^\bullet\circ\alpha^\bullet$
and $\beta^\bullet\circ\phi^\bullet\sim_k\beta^\bullet\circ\psi^\bullet$.
Moreover, $\sim_k$ is an equivalence relation on
$\Hom_{\sFil.\sC(\cA)}(K^\bullet,L^\bullet)$. It follows that there
exists a well defined {\em filtered homotopy category of order $k$}
$$
\sFil.\Hot(\cA,k)
$$
whose objects are the same as those of $\sFil.\sC(\cA)$,
and whose morphisms are the order $k$ homotopy classes of
morphisms of complexes. Furthermore, we have a natural functor
$$
\sFil.\sC(\cA)\to\sFil.\Hot(\cA,k)
$$
which is the identity on objects, and the quotient map on
$\Hom$-sets.
\end{enumerate}
\end{definition}

\sset\subsubsection{}\label{subsec_hot-and-SpSeq}
With the notation of definition \ref{def_fil-homotopy},
suppose that $k\geq 0$, and consider any element
$s^\bullet\in\Fil^k\Hom_\cA^\bullet(K^\bullet,L^\bullet)$; then
$\sigma:=d^{-1}_{K,L}(s^\bullet):
(K^\bullet,\Fil^\bullet K^\bullet)\to
(L^\bullet,\Fil^\bullet L^\bullet)$ is a morphism of
filtered complexes, and a simple inspection shows that
$$
\sigma^{p+q}(Z(K)^{pq}_r)\subset\Ker\,d^{p+q}_{\Fil^pL}
\qquad
\text{for every $p,q,r\in\Z$ with $r\geq k$}
$$
so the induced map $E(\sigma)^{pq}_r:E(K)^{pq}_r\to E(L)^{pq}_r$
vanishes for every $p,q,r\in\Z$ with $r\geq k$. Consequently,
we obtain a commutative diagram of functors
$$
{\diagram
\sFil.\sC(\cA) \ar[r] \ar[d] & \SpSeq_0^\infty(\cA) \ar[d] \\
\sFil.\Hot(\cA,k) \ar[r] & \SpSeq_k^\infty(\cA)
\enddiagram}
\qquad
\text{for every $k\geq 0$}
$$
whose top horizontal arrow is the functor of
\eqref{subsec_Fil-toSpSeq}, and whoses left (resp.
right) vertical arrow is the functor of definition
\ref{def_fil-homotopy}(ii) (resp. the forgetful
functor of \eqref{subsec_simply-forget}).

\sset\subsubsection{}\label{subsec_decalage}
We show next, how to lift the d\'ecalage functor
of \eqref{subsec_decalage-spseq}, to a functor on
complexes
$$
\Dec:\sFil.\sC(\cA)\to\sFil.\sC(\cA).
$$
Namely, given $(K^\bullet,\Fil^\bullet K^\bullet)$
as in \eqref{subsec_fil-spec-seq}, we set
$$
D^\bullet:=K^\bullet
\qquad\text{and}\qquad
\Fil^pD^n:=Z(K)_1^{p+n,-p}
\quad
\text{for every $p,n\in\Z$}.
$$
It follows straightforwardly from \eqref{eq_beg-to-differ},
that the differential $d^\bullet_K$ of $K^\bullet$
restricts to a morphism $\Fil^pD^n\to\Fil^pD^{n+1}$
for every $p,n\in\Z$, so $(D^\bullet,\Fil^\bullet D^\bullet)$
is an object of $\sFil.\sC(\cA)$, and we let
$$
\Dec(K^\bullet,\Fil^\bullet K^\bullet):=
(D^\bullet,\Fil^\bullet D^\bullet).
$$
Every morphism \eqref{eq_fil-map} induces in the obvious
fashion a morphism
$$
\Dec(K^\bullet,\Fil^\bullet K^\bullet)\to
\Dec(L^\bullet,\Fil^\bullet L^\bullet)
$$
so $\Dec$ is a functor as sought. Moreover, it is
easily seen that if $\Fil^\bullet K^\bullet$ is a
bounded filtration, the same holds for the filtration
of $\Dec(K^\bullet,\Fil^\bullet K^\bullet)$.

\sset\subsubsection{}
We wish next to compare the spectral sequences
$$
E(K)^{\bullet\bullet}_\bullet
\qquad\text{and}\qquad
E(\Dec K)^{\bullet\bullet}_\bullet
$$
attached, as in \eqref{subsec_fil-spec-seq}, to
$(K^\bullet,\Fil^\bullet)$ and respectively to
$\Dec(K^\bullet,\Fil^\bullet K)$. To this aim, notice that
$$
Z(\Dec K)_0^{p,n-p}=Z(K)_1^{p+n,-p}
\qquad
D(\Dec K)_0^{p,n-p}=Z(K)_1^{p+1+n,-p-1}
$$
and we have
$$
Z(K)_1^{p+1+n,-p-1}\subset(\Fil^{p+1+n}K^\bullet)^n\subset
D(K)_1^{p+n,-p}\subset Z(K)_1^{p+n,-p}
$$
whence a natural epimorphism
$$
u^{p,n-p}:E(\Dec K)_0^{p,n-p}\to\Dec(E(K))_0^{p,n-p}
\qquad
\text{for every $p,n\in\Z$}.
$$
With this notation, we have :

\begin{proposition}\label{prop_decalage}
With the notation of \eqref{subsec_decalage}, the
following holds :
\begin{enumerate}
\item
The system $(u^{p,n-p}~|~p,n\in\Z)$ extends to an epimorphism
in $\SpSeq_0(\cA)$
$$
u^{\bullet\bullet}_\bullet:E(\Dec K)_\bullet^{\bullet\bullet}\to
\Dec(E(K))_\bullet^{\bullet\bullet}.
$$
\item
The morphism $u^\bullet$ induces isomorphisms in $\SpSeq_1(\cA)$
$$
(E(\Dec K)_r^{\bullet\bullet}~|~r\geq 1)\isom
(\Dec(E(K))_r^{\bullet\bullet}~|~r\geq 1).
$$
\end{enumerate}
\end{proposition}
\begin{proof} For any filtered complex
$(L^\bullet,\Fil^\bullet L^\bullet)$ of $\cA$, let
$$
L'{}^\bullet:=L^\bullet
\qquad
(\Fil^pL'{}^\bullet)^n:=(\Fil^{p-n}L^\bullet)^n
\qquad
\text{for every $p,n\in\Z$}.
$$
It is easily seen that
$$
(L^\bullet,\Fil^\bullet L^\bullet)':=
(L'{}^\bullet,\Fil^\bullet L'{}^\bullet)
$$
is a filtered complex of $\cA$, and a direct inspection
of the definitions yields the identities
$$
E(L)^{p,n-p}_r=E(L')_{r+1}^{p+n,-p}
\qquad
d(L)_r^{p,n-p}=d(L')_{r+1}^{p+n,-p}
\qquad
\text{for every $p,n\in\Z$ and $r\in\N$}
$$
which add up to an identity in $\SpSeq_0(\cA)$ :
\set\begin{equation}\label{eq_filtered-identity}
(E(L)^{\bullet\bullet}_r~|~r\geq 0)=
\Dec(E(L')^{\bullet\bullet}_r~|~r\geq 1).
\end{equation}
Especially, set $(M^\bullet,\Fil^\bullet M^\bullet):=
(\Dec(K^\bullet,\Fil^\bullet K^\bullet))'$. Explicitly,
we have
$$
M^\bullet=K^\bullet
\qquad
(\Fil^pM^\bullet)^n=Z(K)_1^{p,n-p}\subset(\Fil^pK^\bullet)^n
\qquad
\text{for every $p,n\in\Z$}
$$
so that the identity morphism $\one_K:M^\bullet\to K^\bullet$
is a morphism
$$
(M^\bullet,\Fil^\bullet M^\bullet)\to
(K^\bullet,\Fil^\bullet K^\bullet)
\qquad
\text{in $\sFil.\sC(\cA)$}
$$
which in turns yields a morphism in $\SpSeq_0(\cA)$ :
$$
E(M)_\bullet^{\bullet\bullet}\to E(K)^{\bullet\bullet}_\bullet.
$$
Combining with \eqref{eq_filtered-identity}, we get
a natural morphism
$$
v^{\bullet\bullet}_\bullet:E(\Dec K)^{\bullet\bullet}_\bullet\to
\Dec(E(K))^{\bullet\bullet}_\bullet
$$
and a direct inspection shows that
$u^{\bullet\bullet}_\bullet=v^{\bullet\bullet}_\bullet$,
whence (i). To check (ii), we remark that
$$
Z(M)_r^{pq}=
\Ker\,(\bar d:Z(K)^{pq}_1\to Z(K)_1^{p,q+1}/Z(K)_1^{p+r,q-r+1})=
Z(K)^{pq}_r
\qquad
\text{for every $r\geq 1$}
$$
whence :
$$
B(M)^{pq}_r=B(K)^{pq}_r
\qquad
\text{for every $r\geq 1$}
$$
due to lemma \eqref{lem_notice-that}.
\end{proof}

\subsection{Derived categories and derived functors}
\label{sec_derived-cats}
Let $\cA$ be any additive category. With the notation
of remark \ref{rem_chain-homotopies}(ii), set
$$
\C^\bullet:=\Coker\,(\iota^\bullet_1:\Z[0]\to\K\La 1\Ra^\bullet).
$$
Explicitly, $\C^\bullet$ is the object of $\sC^{[-1,0]}(\Z\Mod)$
such that $\C^{-1}=\C^0=\Z$ and with $d_{\C}^{-1}=\one_\Z$.
The morphism $\iota^\bullet_0$ induces a short exact sequence of
complexes of abelian groups :
$$
0\to\Z[0]\xrightarrow{\ \iota^\bullet\ }\C^\bullet
\xrightarrow{\ \pi^\bullet\ }\Z[1]\to 0.
$$
Now, let $\phi^\bullet:K^\bullet\to L^\bullet$ be
any morphism in $\sC(\cA)$. We define the {\em cone}
of $\phi^\bullet$ as the push-out in the cocartesian
diagram of $\sC(\cA)$ :
$$
\xymatrix{
K^\bullet \ar[rr]^-{\iota^\bullet\otimes_\Z\one_K}
\ar[d]_{\phi^\bullet} & &
\C^\bullet\otimes_\Z K^\bullet \ar[d]^{\beta^\bullet} \\
L^\bullet \ar[rr]^-{\psi^\bullet} & & (\Cone\,\phi)^\bullet.
}$$
Let also $\gamma^\bullet:\K\La 1\Ra^\bullet\to\C^\bullet$ be the
natural projection; by remark \ref{rem_chain-homotopies}(ii),
the composition
$$
\beta^\bullet\circ(\gamma^\bullet\otimes_\Z K^\bullet):
\K\La 1\Ra^\bullet\otimes_\Z K^\bullet\to(\Cone\,\phi)^\bullet
$$
corresponds to a homotopy from $\psi^\bullet\circ\phi^\bullet$
to $\beta^\bullet\circ(\gamma^\bullet\circ\iota_1^\bullet)
\otimes_\Z K^\bullet$, and the latter is obviously the
zero endomorphism of $K^\bullet$. Moreover, we get a
natural identification
\set\begin{equation}\label{eq_coker-is-shift}
\Coker\,\psi^\bullet\isom
\Coker\,(\iota^\bullet\otimes_\Z K^\bullet)\isom K^\bullet[1]
\end{equation}
and we let $\partial^\bullet$ be the unique morphism
of complexes which fits in the commutative diagram
$$
\xymatrix{ & (\Cone\,\phi)^\bullet \ar[ld]_{-\pi^\bullet}
\ar[rd]^{\partial^\bullet} \\
\Coker\,\psi^\bullet \ar[rr] & & K^\bullet[1]
}$$
whose horizontal arrow is the identification
\eqref{eq_coker-is-shift}, and where $\pi^\bullet$
is the natural projection. Summing up, we have attached
to $\phi^\bullet$ a natural sequence of morphisms of complexes
$$
\Theta(\phi^\bullet)
\quad :\quad
K^\bullet\xrightarrow{\ \phi^\bullet\ }L^\bullet
\xrightarrow{\ \psi^\bullet\ }(\Cone\,\phi)^\bullet
\xrightarrow{\ \partial^\bullet\ }K^\bullet[1]
$$
such that $\psi^\bullet\circ\phi^\bullet$ is homotopically
trivial, and the induced sequence of morphisms in $\cA$
$$
0\to L^i\xrightarrow{\ \psi^i\ }(\Cone\,\phi)^i
\xrightarrow{\ \partial^i\ } K^{i+1}\to 0
$$
is {\em split exact} for every $i\in\Z$, {\em i.e.}
such that $(\Cone\,\phi)^i=L^i\oplus K^{i+1}$, and
$\psi_i$ (resp. $-\partial^i$) is the natural monomorphism
(resp. epimorphism) induced by this direct sum
decomposition. A direct inspection shows that the
differential
$$
d_{\Cone\,\phi}^{i-1}:
L^{i-1}\oplus K^i\to L^i\oplus K^{i+1}
$$
is given by the matrix
$$
\left[\begin{array}{cc}
d^{i-1}_L & \phi^i \\
0 & -d^i_K
\end{array}\right].
$$
Moreover, it is easily seen that every morphism
$\beta^\bullet:\phi_1^\bullet\to\phi^\bullet_2$ in
$\sMorph(\sC(\cA))$
induces a natural commutative diagram in $\sC(\cA)$ :
\set\begin{equation}\label{eq_cone-is-nature}
{\diagram
\Theta(\phi_1^\bullet) \ar[d]_{\Theta(\beta^\bullet)} &
K_1^\bullet \ar[rr]^-{\phi_1^\bullet} \ar[d]_{\beta_1^\bullet} & &
L_1^\bullet \ar[rr]^-{\psi_1^\bullet} \ar[d]_{\beta_2^\bullet} & &
(\Cone\,\phi_1)^\bullet \ar[d]_{\Cone(\beta_1,\beta_2)^\bullet}
\ar[rr]^-{\partial_1^\bullet} & &
K^\bullet_1[1] \ar[d]^{\beta_1^\bullet[1]} \\
\Theta(\phi^\bullet_2) &
K_2^\bullet \ar[rr]^-{\phi_2^\bullet} & &
L_2^\bullet \ar[rr]^-{\psi_2^\bullet} & &
\ar[rr]^-{\partial_2^\bullet} (\Cone\,\phi_2)^\bullet & &
K^\bullet_2[1]
\enddiagram}
\end{equation}
where $\psi^\bullet_1$ and $\psi^\bullet_2$ are the
natural morphisms, as in the foregoing (details left
to the reader).

\begin{definition}\label{def_il-triangolo-no}
Let $\cA$ be any additive category, and denote by
$\sT(\cA)$ either of the categories $\sC(\cA)$ or
$\Hot(\cA)$. With the notation of \eqref{sec_derived-cats} :
\begin{enumerate}
\item
We call $\Theta(\phi^\bullet)$ the {\em true triangle}
associated to the morphism $\phi^\bullet$, and
$\partial^\bullet$ is the {\em boundary morphism}
of $\Theta(\phi^\bullet)$.
\item
A {\em triangle} of $\sT(\cA)$ is any sequence of morphisms
of $\sT(\cA)$ :
\set\begin{equation}\label{eq_renato-zero}
A^\bullet\xrightarrow{\ \alpha^\bullet\ }
B^\bullet\xrightarrow{\ \beta^\bullet\ }
C^\bullet\xrightarrow{\ \gamma^\bullet\ }
A^\bullet[1].
\end{equation}
\item
Let
$\Theta_i:=(A^\bullet_i\xrightarrow{\ \alpha^\bullet_i\ }
B^\bullet_i\xrightarrow{\ \beta^\bullet_i\ }
C^\bullet_i\xrightarrow{\ \gamma^\bullet_i\ }A^\bullet_i[1])$
for $i=1,2$ be two triangles of $\sT(\cA)$.
A {\em morphism of triangles} $\Theta_1\to\Theta_2$
is a commutative diagram of morphisms of $\sT(\cA)$ :
$$
\xymatrix{ 
A^\bullet_1 \ar[r]^-{\alpha_1^\bullet} \ar[d]_{\tau^\bullet} &
B^\bullet_1 \ar[r]^-{\beta_1^\bullet} \ar[d] &
C^\bullet_1 \ar[r]^-{\gamma_1^\bullet} \ar[d] &
A^\bullet_1[1] \ar[d]^{\tau^\bullet[1]} \\
A^\bullet_2 \ar[r]^-{\alpha_2^\bullet} &
B^\bullet_2 \ar[r]^-{\beta_2^\bullet} &
C^\bullet_2 \ar[r]^-{\gamma_2^\bullet} &
A^\bullet_2[1].
}$$
Morphisms of triangles can be composed in the obvious
fashion, and clearly the system of all triangles of
$\sT(\cA)$ and their morphisms, forms a category.
\item
A {\em distinguished triangle} of $\sT(\cA)$ is a
triangle of $\sT(\cA)$ that is isomorphic to a true
triangle. We denote by
$$
\Theta.\sT(\cA)
$$
the full subcategory of the category of triangles of
$\sT(\cA)$ whose objects are the distinguished triangles.
\item
Let $\cB$ be any other additive category. A {\em triangulated
functor} from $\sT(\cA)$ to $\sT(\cB)$ is a pair $(F,\tau)$
consisting of a functor $F:\sT(\cA)\to\sT(\cB)$ and a
natural isomorphism
$$
\tau^\bullet_K:F(K^\bullet[1])\isom(FK^\bullet)[1]
\qquad
\text{for every $K^\bullet\in\Ob(\sT(\cA))$}
$$
such that, for every distinguished triangle
\eqref{eq_renato-zero} of $\sT(\cA)$, the triangle
$$
FA^\bullet\xrightarrow{\ F\alpha^\bullet\ }FB^\bullet
\xrightarrow{\ F\beta^\bullet\ }FC^\bullet
\xrightarrow{\ \tau^\bullet_A\circ F\gamma^\bullet\ }
FA^\bullet[1]
$$
is distinguished in $\sT(\cB)$. We also say that $F$ is
{\em triangulated} if there exists a natural isomorphism
$\tau$ as above, such that $(F,\tau)$ is a triangulated
functor.
\end{enumerate}
\end{definition}

\begin{remark}\label{rem_long-exact-triang}
(i)\ \
With the terminology of definition \ref{def_il-triangolo-no},
we may say that diagram \eqref{eq_cone-is-nature} is
a morphism of distinguished triangles, and clearly
the rules : $\phi^\bullet\mapsto\Theta(\phi^\bullet)$
for any morphism $\phi^\bullet$ in $\sC(\cA)$ and
$\beta^\bullet\mapsto\Theta(\beta^\bullet)$ for any
morphism $\beta^\bullet$ in $\sMorph(\sC(\cA))$
amount to a functor
\set\begin{equation}\label{eq_can-improve-this}
\sMorph(\sC(\cA))\to\Theta.\sC(\cA).
\end{equation}
Moreover, if $\phi^\bullet$ is a morphism of
$\sC^+(\cA)$ (resp. $\sC^-(\cA)$, resp. $\sC^b(\cA)$),
then clearly $(\Cone\,\phi)^\bullet$ lies in the same
subcategory of $\sC(\cA)$, so \eqref{eq_can-improve-this}
restricts to functors
$$
\sMorph(\sC^?(\cA))\to\Theta.\sC^?(\cA)
\qquad
\text{with ``$?$'' equal to either $+$, $-$, or $b$}
$$
with obvious notation.

(ii)\ \
Let $F:\cA\to\cB$ be any additive functor between additive
categories. Directly from the definitions, it is clear that
the induced functors $\sC(F)$ and $\Hot(F)$ are both
triangulated.

(iii)\ \
Suppose now that $\cA$ is an abelian category, so that
the same holds for $\sC(\cA)$, and the discussion of
\eqref{sec_derived-cats} yields a short exact sequence
$$
0\to L^\bullet\xrightarrow{\ \psi^\bullet\ }
(\Cone\,\phi^\bullet)\xrightarrow{\ \partial^\bullet\ }
K^\bullet[1]\to 0
\qquad
\text{in $\sC(\cA)$}
$$
such that the boundary map in degree $i$ of the induced
long exact cohomology sequence is none else than
$\phi^{i+1}$, so we get a natural acyclic complex 
$$
\cdots\to H^{i-1}L^\bullet\xrightarrow{\ \psi^i\ }
H^{i-1}(\Cone\,\phi)^\bullet
\xrightarrow{\ \partial^{i-1}\ }H^iK^\bullet
\xrightarrow{\ \phi^i\ }H^iL^\bullet\to\cdots
$$
for every morphism $\phi^\bullet$ in $\sC(\cA)$.
However, one of the subtleties of distinguished triangles,
is that they provide a language for expressing certain
exactness and cohomological assertions, that remains
available even for general additive but not necessarily
abelian categories. For instance, let us show the following :
\end{remark}

\begin{proposition}\label{prop_triangul-hot}
Let $\cA$ be any additive category. For any complex
$L^\bullet$ of $\cA$, the functors
$$
\begin{aligned}
\sC(\cA)\to\, & \sC(\Z\Mod)
& \qquad &
K^\bullet\mapsto\Hom^\bullet_\cA(L^\bullet,K^\bullet)
\\
\Hot(\cA)\to\, & \Hot(\Z\Mod)
& \qquad &
K^\bullet\mapsto\rHot^\bullet_\cA(L^\bullet,K^\bullet)
\end{aligned}
$$
are triangulated (notation of example
{\em\ref{ex_hot-morphisms}}).
\end{proposition}
\begin{proof} Let $\phi^\bullet:K^\bullet_1\to K^\bullet_2$ be
any morphism of complexes of $\cA$. The universal properties
of the push-out and of the cokernel yield a unique commutative
diagram
$$
\xymatrix@C-5pt{
\Hom^\bullet_\cA(L^\bullet,\Theta(\phi^\bullet)) \ar[d] &
\Hom^\bullet_\cA(L^\bullet,K_2^\bullet) \ar[r] \ddouble &
\Hom^\bullet_\cA(L^\bullet,(\Cone\,\phi^\bullet))
\ar[r] \ar[d]^{\beta^\bullet} &
\Hom^\bullet_\cA(L^\bullet,K_1^\bullet[1]) \ar[d]^{\gamma^\bullet}
\\
\Theta(\Hom^\bullet_\cA(L^\bullet,\phi^\bullet)) &
\Hom^\bullet_\cA(L^\bullet,K_2^\bullet) \ar[r] &
\Cone\,\Hom^\bullet_\cA(L^\bullet,\phi^\bullet) \ar[r] &
\Hom^\bullet_\cA(L^\bullet,K_1^\bullet)[1].
}$$
Now, since the sequence
$0\to K_2^i\to(\Cone\,\phi)^i\to K_1^{i+1}\to 0$ is split
exact in each degree $i\in\Z$, it is easily seen that
the same holds for the induced sequence
$$
0\to \Hom^i_\cA(L^\bullet,K_2^\bullet)\to
\Hom^i_\cA(L^\bullet,(\Cone\,\phi)^\bullet)
\to\Hom^i_\cA(L^\bullet,K^\bullet_1[1])\to 0.
$$
It follows that $\gamma^\bullet$ is the obvious natural
identification and $\beta^\bullet$ is the direct product
of the natural identifications
$$
\Hom_\cA(L^i,\phi^j\oplus\phi^{j+1})\isom
\Hom_\cA(L^i,\phi^j)\oplus\Hom_\cA(L^i,\phi^{j+1})
\qquad
\text{for every $i,j\in\Z$}
$$
whence the assertion for the functor $\Hom^\bullet_\cA$.
The assertion for $\rHot^\bullet_\cA$ follows immediately.
\end{proof}

\begin{remark}\label{rem_biadditive-triang}
Resume the situation of \eqref{subsec_hot-biadditive}.
It is easily seen that, {\em mutatis mutandis}, the
proof of proposition \ref{prop_triangul-hot} shows also
the following. For every $K^\bullet\in\Ob(\Hot(\cA))$
the functor
$$
\Hot^-(\cA')\to\Hot^-(\cA'')
\qquad
L^\bullet\mapsto B^\bullet_-(K^\bullet,L^\bullet)
$$
is triangulated. Likewise, the same holds for the functor
$$
\Hot^-(\cA)\to\Hot^-(\cA'')
\qquad
K^\bullet\mapsto B^\bullet_-(K^\bullet,L^\bullet)
$$
for every $L^\bullet\in\Ob(\Hot(\cA'))$ : details left
to the reader.
\end{remark}

\begin{lemma}\label{lem_first-cone-lemma}
With the notation of \eqref{sec_derived-cats}, let
$\phi^\bullet_1,\phi^\bullet_2:K^\bullet\to L^\bullet$
be any two morphisms in $\sC(\cA)$. The following holds :
\begin{enumerate}
\item
Any homotopy $s^\bullet$ from $\phi_1^\bullet$ to
$\phi_2^\bullet$ induces a natural isomorphism
$$
\gamma^\bullet_s:
(\Cone\,\phi_1)^\bullet\isom(\Cone\,\phi_2)^\bullet
\qquad
\text{in $\sC(\cA)$}
$$
fitting into a commutative diagram
$$
\xymatrix{
K^\bullet \ar[r]^-{\phi_1^\bullet} &
L^\bullet \ar[rr]^-{\psi_1^\bullet} \ddouble & &
(\Cone\,\phi_1)^\bullet
\ar[rr]^-{\partial^\bullet_{\phi_1}} \ar[d]_{\gamma_s^\bullet} & &
K^\bullet[1] \ddouble \\
K^\bullet \ar[r]^-{\phi_2^\bullet} &
L^\bullet \ar[rr]^-{\psi_2^\bullet} & &
(\Cone\,\phi_2)^\bullet \ar[rr]^-{\partial^\bullet_{\phi_2}} & &
K^\bullet[1]
}$$
whose top (resp. bottom) row is $\Theta(\phi^\bullet_1)$ (resp.
$\Theta(\phi^\bullet_2)$).
\item
There is a natural isomorphism
$$
\omega^\bullet:(\Cone\,\psi_1)^\bullet\isom
K^\bullet[1]\oplus(\C^\bullet\otimes_\Z L^\bullet)
\qquad
\text{in $\sC(\cA)$}
$$
fitting into a commutative diagram
$$
\xymatrix{ L^\bullet \ar[r]^-{\psi_1^\bullet} &
(\Cone\,\phi_1)^\bullet
\ar[r]^-{\tau^\bullet} \ar[d]_{\partial^\bullet_{\phi_1}} &
(\Cone\,\psi_1)^\bullet
\ar[d]^{\omega^\bullet} \ar[r]^-{\partial^\bullet_{\psi_1}} &
L^\bullet[1] \ddouble \\
& K^\bullet[1] & \ar[l]_-{-p_1^\bullet}
K^\bullet[1]\oplus(\C^\bullet\otimes_\Z L^\bullet)
\ar[r]^-{p^\bullet_2} & L^\bullet[1]
}$$
whose top row is $\Theta(\psi_1^\bullet)$, and
where $p_1^\bullet$ is the natural projection,
and $p_2^\bullet$ is the morphism given by the matrix
$$
\left[\begin{array}{c}
0 \\
\phi^{i+1}_1 \\
-\one_{L^{i+1}}
\end{array}\right]:
L^i\oplus K^{i+1}\oplus L^{i+1}\to L^{i+1}
\qquad
\text{for every $i\in\Z$}.
$$
\end{enumerate}
\end{lemma}
\begin{proof}(i): By assumption,
$\phi_1^i=\phi_2^i+d^{i-1}s^i+s^{i+1}d^i$ for every $i\in\Z$.
We let $\gamma^i_s$ be the automorphism of $L^i\oplus K^{i+1}$
given by the matrix
$$
\left[\begin{array}{cc}
\one_{L^i} & s^{i+1} \\
0 & \one_{K^{i+1}}
\end{array}\right]
\qquad
\text{for every $i\in\Z$}.
$$
A direct computation shows that the system
$(\gamma^i_s~|~i\in\Z)$ is the sought isomorphism, and
the commutativity of the resulting diagram as in (i)
follows by a simple inspection.

(ii): Set $\sigma^\bullet:=
(\iota^\bullet\otimes_\Z\one_L)\circ\phi^\bullet_1$.
By definition, we have a natural identification
$$
(\Cone\,\sigma)^\bullet=(\Cone\,\psi_1)^\bullet
$$
as well as a commutative diagram
$$
\xymatrix{
K^\bullet \ar[r]^-{\phi_1^\bullet} \ar[d]_{\iota^\bullet\otimes_\Z\one_K}
& L^\bullet \ar[r]^-{\iota^\bullet\otimes_\Z\one_L}
\ar[d]_{\psi_1^\bullet} &
\C^\bullet\otimes_\Z L^\bullet \ar[d] \\
\C^\bullet\otimes_\Z K^\bullet \ar[r] \ar[d] &
(\Cone\,\phi_1)^\bullet
\ar[r]^-{\tau^\bullet} \ar[d]_{\partial^\bullet_{\phi_1}} &
(\Cone\,\psi_1)^\bullet \ar[d]_{\partial^\bullet_\sigma}
\\
K^\bullet[1] \ar[r]^-{-\one_{K[1]}} & K^\bullet[1] \rdouble
& K^\bullet[1]
}$$
whose two square subdiagrams on the top ladder are cocartesian.
Notice that $\C^\bullet$ is homotopically trivial (details
left to the reader); by remark \ref{rem_chain-homotopies}(iv),
the same then holds for $\C^\bullet\otimes_\Z L^\bullet$.
Thus, $\sigma^\bullet$ is null-homotopic.
In light of (i) and \eqref{eq_coker-is-shift}, we deduce an
isomorphism as sought, fitting into a commutative diagram
$$
\xymatrix{ K^\bullet \ar[r] &
\C^\bullet\otimes_\Z L^\bullet \ar[r] \ddouble &
(\Cone\,\psi_1)^\bullet \ar[d]^{\omega^\bullet}
\ar[r]^-{\partial^\bullet_\sigma} & K^\bullet[1] \ddouble
\\
K^\bullet \ar[r]^-0 & \C^\bullet\otimes_\Z L^\bullet \ar[r] &
K^\bullet[1]\oplus(\C^\bullet\otimes_\Z L^\bullet)
\ar[r]^-{\partial_0^\bullet} &  K^\bullet[1]
}$$
whose bottom row is the true triangle associated to the
zero morphism $K^\bullet\to\C^\bullet\otimes_\Z L^\bullet$.
Explicitly, the differentials
$$
d_{(\Cone\,\psi_1)}^i,
d_{K[1]\oplus(\C\otimes_\Z L)}^i:
L^i\oplus K^{i+1}\oplus L^{i+1}\to
L^{i+1}\oplus K^{i+2}\oplus L^{i+2}
$$
are given by the matrices
$$
\left[\begin{array}{ccc}
d_L^i & \phi^{i+1}_1 & \one_{L^{i+1}} \\
0 & -d_K^{i+1} & 0 \\
0 & 0 & -d_L^{i+1}
\end{array}\right]
\qquad
\left[\begin{array}{ccc}
d_L^i & 0 & \one_{L^{i+1}} \\
0 & -d_K^{i+1} & 0 \\
0 & 0 & -d_L^{i+1}
\end{array}\right]
\qquad
\text{for every $i\in\Z$}
$$
from which we see that we can take for $\omega^i$ the
automorphism of $L^i\oplus K^{i+1}\oplus L^{i+1}$ given
by the matrix
$$
\left[\begin{array}{ccc}
\one_{L^i} & 0 & 0 \\
0 & \one_{K^{i+1}} & 0 \\
0 & \phi^{i+1}_1 & \one_{L^{i+1}}
\end{array}\right]
\qquad
\text{for every $i\in\Z$}
$$
and then the commutativity of the right square
subdiagram of (ii) is immediate. Lastly, it is easily
seen that $\partial_0^\bullet=-p^\bullet_1$, therefore
$$
-p^\bullet_1\circ\omega^\bullet\circ\tau^\bullet=
\partial^\bullet_\sigma\circ\tau^\bullet=
\partial^\bullet_{\phi_1}
$$
which shows the commutativity of the left square
subdiagram of (ii).
\end{proof}

\begin{remark}\label{rem_hom-triangles}
(i)\ \
In the situation of lemma \ref{lem_first-cone-lemma},
notice that the diagram
$$
\xymatrix{
K^\bullet \ar[r]^-{\phi_1^\bullet} \ddouble &
L^\bullet \ar[rr]^-{\psi_1^\bullet} \ddouble & &
(\Cone\,\phi_1)^\bullet
\ar[rr]^-{\partial^\bullet_{\phi_1}} \ar[d]_{\gamma_s^\bullet} & &
K^\bullet[1] \ddouble \\
K^\bullet \ar[r]^-{\phi_2^\bullet} &
L^\bullet \ar[rr]^-{\psi_2^\bullet} & &
(\Cone\,\phi_2)^\bullet \ar[rr]^-{\partial^\bullet_{\phi_2}} & &
K^\bullet[1]}
$$
commutes in $\Hot(\cA)$, and amounts to a natural
isomorphism
$$
\Theta(\phi_1^\bullet)\isom\Theta(\phi_2^\bullet)
\qquad
\text{in $\Theta.\Hot(\cA)$}.
$$

(ii)\ \
Let $\Theta:=(A^\bullet\xrightarrow{\ \alpha^\bullet\ }
B^\bullet\xrightarrow{\ \beta^\bullet\ }
C^\bullet\xrightarrow{\ \gamma^\bullet\ }A^\bullet[1])$
be any triangle of $\sC(\cA)$ or $\Hot(\cA)$. We obtain
a new triangle by setting
$$
\Theta[1]:=(B^\bullet\xrightarrow{\ \beta^\bullet\ }
C^\bullet\xrightarrow{\ \gamma^\bullet\ }A^\bullet[1]
\xrightarrow{\ -\alpha^\bullet[1]\ }B^\bullet[1])
$$
and clearly, the rule $\Theta\mapsto\Theta[1]$ yields
an endofunctor of the category of triangles.
The change of sign on the last arrow is needed to ensure
that this operator restricts to an endofunctor of the
subcategory of distinguished triangles, at least 
up to homotopy. Indeed, we may state :
\end{remark}

\begin{proposition}\label{prop_rotate-triangles}
With the notation of remark {\em\ref{rem_hom-triangles}(ii)},
suppose that $\Theta$ is a distinguished triangle of\/
$\Hot(\cA)$. Then the same holds for $\Theta[1]$, so
we have an automorphism
$$
\Theta.\Hot(\cA)\isom\Theta.\Hot(\cA)
\qquad
\Theta\mapsto\Theta[1].
$$
\end{proposition}
\begin{proof} We may assume that $\Theta=\Theta(\alpha^\bullet)$
for some morphism $\alpha^\bullet:A^\bullet\to B^\bullet$ of
complexes of $\cA$, in which case $C^\bullet=(\Cone\,\alpha)^\bullet$,
and it suffices to show the more precise :

\begin{claim} There is a natural isomorphism
$$
\Theta(\beta^\bullet)\isom\Theta(\alpha^\bullet)[1]
\qquad
\text{in $\Theta.\Hot(\cA)$}.
$$
\end{claim}
\begin{pfclaim}[] Indeed, notice that the natural projection
$p_1^\bullet:A^\bullet[1]\oplus(\C^\bullet\otimes_\Z B^\bullet)
\to A^\bullet[1]$ represents an isomorphism in $\Hot(\cA)$,
whose inverse is the class of the natural monomorphism
$e_1^\bullet: A^\bullet[1]\to
A^\bullet[1]\oplus(\C^\bullet\otimes_\Z B^\bullet)$. Taking
into account lemma \ref{lem_first-cone-lemma}(ii), we then
get the commutative diagram
in $\Hot(\cA)$ :
$$
\xymatrix{ B^\bullet \ar[r]^-{\beta^\bullet} \ddouble &
(\Cone\,\alpha)^\bullet \ddouble \ar[r] &
(\Cone\,\beta)^\bullet \ar[d]^{p_1^\bullet\circ\omega^\bullet}
\ar[r]^-{\partial_\beta^\bullet} &
B^\bullet[1] \ddouble
\\
B^\bullet \ar[r]^-{\beta^\bullet} &
(\Cone\,\alpha)^\bullet \ar[r]^-{\partial_\alpha^\bullet} &
A^\bullet[1] \ar[r]^-{-\alpha^\bullet} & B^\bullet[1]
}$$
where $\omega^\bullet:(\Cone\,\beta)^\bullet\isom
A^\bullet[1]\oplus(\C^\bullet\otimes_\Z B^\bullet)$
is the isomorphism provided by lemma
\ref{lem_first-cone-lemma}(ii). The claim follows.
\end{pfclaim}
\end{proof}

\begin{definition}\label{def_derived-cat}
Let $\cA$ be any abelian category.

(i)\ \
A {\em quasi-isomorphism} is a morphism
$\phi^\bullet:K^\bullet\to L^\bullet$ in $\sC(\cA)$
such that $H^i\phi^\bullet$ is an isomorphism for
every $i\in\Z$.

(ii)\ \
With the notation of \eqref{sec_derived-cats}, we let
$\Sigma$ be the subset of $\mathrm{Morph}(\sC(\cA))$
consisting of all quasi-isomorphisms. We set 
$$
\sD(\cA):=\sC(\cA)[\Sigma^{-1}]
$$
(see remark \ref{rem_local-cat}(i)) and we call this category
the {\em derived category} of $\cA$.

(iii)\ \
More generally, if $I\subset\Z$ is any interval as in
\eqref{sec_brutal-truncate}, and $a\in\Z$ is any integer,
we denote
$$
\sD^I(\cA)
\qquad
\sD^{\geq a}(\cA)
\qquad
\sD^{\leq a}(\cA)
\qquad
\sD^+(\cA)
\qquad
\sD^-(\cA)
\qquad
\sD^b(\cA)
$$
the essential image in $\sD(\cA)$ of the categories
$\sC^I(\cA)$, respectively $\sC^{\geq a}(\cA)$, 
$\sC^{\leq a}(\cA)$, $\sC^+(\cA)$, $\sC^-(\cA)$,
$\sC^b(\cA)$.
\end{definition}

\begin{remark}\label{rem_derived-cat}
(i)\ \
Notice that, if $\phi^\bullet$ is a quasi-isomorphism
in $\sC(\cA)$, and $\psi^\bullet$ is any other morphism
that is homotopy equivalent to $\phi^\bullet$, then
$\psi^\bullet$ is a quasi-isomorphism as well (remark
\ref{rem_abel-homotopies}(ii)). We may then say that
a morphism $\phi^\bullet$ in $\Hot(\cA)$ is a
{\em quasi-isomorphism}, if it admits a representative
in $\sC(\cA)$ which is a quasi-isomorphism; by the
foregoing, this property can be checked on any
representative for the homotopy class $\phi^\bullet$.

(ii)\ \
In the same vein, notice that, by virtue of proposition
\ref{prop_univ-hot}, the localization functor
$\sC^?(\cA)\to\sD^?(\cA)$ factors uniquely through a functor
\set\begin{equation}\label{eq_hot-abel}
\Hot^?(\cA)\to\sD^?(\cA)
\end{equation}
whenever ``$?$'' equals either $+$, $-$, $b$, or any interval
$I\subset\Z$.

(iii)\ \
In light of remark \ref{rem_abel-homotopies}(iii), we see
that the isomorphism of remark \ref{rem_homotopies}(iii)
identifies the subset $\Sigma^o\subset\Ob(\sC(\cA)^o)$
with the system of quasi-isomorphisms of $\sC(\cA^o)$.
Combining with remark \ref{rem_local-cat}(ii) we get
natural isomorphisms of categories :
$$
\sD^I(\cA)^o\isom\sD^{-I}(\cA^o)
\qquad
\sD^+(\cA)^o\isom\sD^-(\cA^o)
\qquad
\sD^b(\cA)^o\isom\sD^b(\cA^o)
$$
where $I\subset\Z$ is any interval, and we set
$-I:=\{-a~|~a\in I\}$.

(iv)\ \
Notice that the brutal truncation functors of
\eqref{subsec_brutal-truncate} do not transform
quasi-isomorphism into quasi-isomorphisms, hence
they do not descend to the derived category. On
the other hand, the normalized truncation functors
do preserve quasi-isomorphisms, so they yield
functors
$$
\tau^{\geq a}:\sD(\cA)\to\sD^{\geq a}(\cA)
\qquad\text{and}\qquad
\tau^{\leq a}:\sD(\cA)\to\sD^{\leq a}(\cA).
\qquad
\text{for every $a\in\Z$}
$$
Moreover, let $i^{\geq a}:\sC^{\geq a}(\cA)\to\sC(\cA)$
and $j^{\geq a}:\sD^{\geq a}(\cA)\to\sD(\cA)$ be the
inclusion functors; so the counit of the adjoint pair
$(\tau^{\geq a},i^{\geq a})$ is just the identity endofunctor
$\one_{\sC^{\geq a}(\cA)}$, and it is easily seen that the unit
of adjunction  $\one_{\sC(\cA)}\to i^{\geq a}\circ\tau^{\geq a}$
induces a natural transformation
$\eta:\one_{\sD(\cA)}\to j^{\geq a}\circ\tau^{\geq a}$.
Then, the triangular identities for the adjunction
$(\tau^{\geq a},i^{\geq a})$ imply corresponding
triangular identities for $\eta$ and the identity
endofunctor $\one_{\sD^{\geq a}(\cA)}$; by proposition
\ref{prop_triangular-identities}(i), we conclude that
$\tau^{\geq a}$ is left adjoint to $j^{\geq a}$. Likewise,
$\tau^{\leq a}$ is right adjoint to the inclusion
functor $\sD^{\leq a}(\cA)\to\sD(\cA)$.

(v)\ \
Directly from the definition (and the universal property
of localization), we see that, for every $i\in\Z$, the
cohomology functor in degree $i$ factors uniquely through
a functor
$$
H^i:\sD(\cA)\to\cA.
$$
\end{remark}

\begin{theorem}\label{th_from-hot-to-der}
Let $\cA$ be any abelian category, and denote by
$$
\Sigma_?\subset\mathrm{Morph}(\Hot^?(\cA))
\qquad
\text{with ``$?$'' equal to either $+$, $-$, $b$, or any
interval $I\subset\Z$}
$$
the subset of all quasi-isomorphisms. We have :
\begin{enumerate}
\item
The functor \eqref{eq_hot-abel} factors uniquely through
an equivalence :
$$
\Hot^?(\cA)[\Sigma_?^{-1}]\isom\sD^?(\cA).
$$
\item
The set $\Sigma_?$ admits both a right and a left calculus
of fractions.
\end{enumerate}
\end{theorem}
\begin{proof} To begin with, we point out :

\begin{claim}\label{cl_question-marks}
For any abelian category $\cA$, the following holds :
\begin{enumerate}
\addenu\addenu
\item
Let $K^\bullet,L^\bullet$ be any two complexes of $\cA$,
with $L^\bullet\in\Ob(\sC^{\leq a}(\cA))$. If there exists
a quasi-isomorphism $K^\bullet\to L^\bullet$, then the
counit of adjunction $\tau^{\leq a}K^\bullet\to K^\bullet$
is a quasi-isomorphism.
\item
In order to prove the theorem, it suffices to show
that $\Sigma$ admits a right calculus of fractions.
\end{enumerate}
\end{claim}
\begin{pfclaim}(iii) follows immediately from remark
\ref{rem_derived-cat}(iv) (details left to the reader).

(iv): Indeed, suppose that $\Sigma$ admits a right
calculus of fraction; in view of remark
\ref{rem_derived-cat}(iii) it follows that $\Sigma$
admits as well a left calculus of fraction.
Next, in case ``$?$'' equals $\Z$, assertion (i) says that
the induced functor $\Hot(\cA)[\Sigma^{-1}]\isom\sD(\cA)$
is an equivalence; the latter follows immediately by
comparing the respective universal properties : details
left to the reader. So, the theorem is completely proven,
for ``$?$'' equal to $\Z$.

Next, from (iii), the case ``$?$''$=\Z$ of the theorem,
and proposition \ref{prop_full-subcat-fract} we get
assertion (i) of the theorem for ``$?$'' equal to
$]-\infty,a]$, and we also deduce that $\Sigma_{]-\infty,a]}$
admits a right calculus of fractions. We check directly
that $\Sigma_{]-\infty,a]}$ admits also a left calculus
of fractions. Indeed, (the duals of) axioms (CF1) and
(CF2) are obviously fulfilled. For (CF3), let us consider
any pair of morphisms $f^\bullet:B^\bullet\to A^\bullet$
and $s^\bullet:B^\bullet\to C^\bullet$ in
$\Hot^{]-\infty,a]}(\cA)$, with $s^\bullet$ a
quasi-isomorphism; since (the dual of) (CF3) is already
known to hold for $\Sigma$, we may find morphisms
$t:A\to D$ and $g:C\to D$ in $\Hot(\cA)$, such that
$t\circ f=g\circ s$, and where $t$ is a quasi-isomorphism.
But then, by remark \ref{rem_derived-cat}(iv), we may
replace $D$ (resp. $t$, resp. $g$) by $\tau^{\leq a}D$
(resp. by $\tau^{\leq a}t$, resp. by $\tau^{\leq a}g$),
after which we may assume that $t\in\Sigma_{]-\infty,a]}$,
whence (CF3). The proof of (CF4) is similar, and shall
be left to the reader.

Dualizing the foregoing case, and taking into account
remark \ref{rem_derived-cat}(iii), we then get both
assertions (i) and (ii) of the theorem also for ``$?$''
equal to $[a,+\infty[$, so the proof of theorem is complete
for the case when ``$?$'' equals any unbounded interval
$I\subset\Z$. One argues likewise to show the theorem in
case ``$?$'' equals either $+$ or $-$ (details left to
the reader).
Lastly, for ``$?$'' equal to either $b$ or an interval
$[a,b]$, it suffices to apply proposition
\ref{prop_full-subcat-fract} and (iii) to the fully
faithful inclusion functors
$$
\Hot^{[a,b]}(\cA)\to\Hot^{\geq a}(\cA)
\qquad
\Hot^b(\cA)\to\Hot^+(\cA)
$$
and argue as in the foregoing cases to complete the
proof of the claim.
\end{pfclaim}

By claim \ref{cl_question-marks}(ii), it remains only
to show that $\Sigma$ admits a right calculus of
fractions. To this aim, we check the conditions of
definition \ref{def_right-calculus}. (CF1) and (CF2)
are obvious. Next, consider two morphisms
$\phi^\bullet:A^\bullet\to B^\bullet$ and
$\psi^\bullet:C^\bullet\to B^\bullet$ in $\sC^?(\cA)$,
with $\phi^\bullet$ a quasi-isomorphism. Let
$\beta^\bullet:B^\bullet\to(\Cone\,\psi)^\bullet$ be
the natural morphism; with the notation of
\eqref{eq_cone-is-nature}, we have an induced commutative
diagram
$$
\xymatrix{ A^\bullet \ar[rr]^-{\beta^\bullet\circ\phi^\bullet}
\ar[d]_{\phi^\bullet} & & (\Cone\,\psi)^\bullet \ar[r] \ddouble &
(\Cone\,\beta\circ\phi)^\bullet
\ar[d]_{\gamma^\bullet} \ar[r] &
A^\bullet[1] \ar[d]^{\phi^\bullet[1]} \\
B^\bullet \ar[rr]^-{\beta^\bullet} & & (\Cone\,\psi)^\bullet \ar[r] &
(\Cone\,\beta)^\bullet \ar[r] & B^\bullet[1]
}$$
where $\gamma^\bullet:=\Cone\,(\phi^\bullet,\one)^\bullet$,
and by the 5-lemma, it follows that $\gamma^\bullet$ is
a quasi-isomorphism. However, lemma \ref{lem_first-cone-lemma}(ii)
identifies $\gamma^\bullet$ with a morphism
$(\Cone\,\beta\circ\phi)^\bullet\to
C^\bullet[1]\oplus(\C^\bullet\otimes_\Z B^\bullet)$, whence
a commutative diagram
$$
\xymatrix{ (\Cone\,\beta\circ\phi)^\bullet
\ar[d]_{\pi^\bullet\circ\gamma^\bullet} \ar[r] &
A^\bullet[1] \ar[d]^{\phi^\bullet[1]} \\
C^\bullet[1] \ar[r]^-{\psi^\bullet[1]} & B^\bullet[1]
}$$
where $\pi^\bullet:
C^\bullet[1]\oplus(\C^\bullet\otimes_\Z B^\bullet)\to C^\bullet[1]$
is the natural projection. Since
$\C^\bullet\otimes_\Z B^\bullet$ is homotopically trivial
(see the proof of lemma \ref{lem_first-cone-lemma}(ii)),
we see that $\pi^\bullet$ is also a quasi-isomorphism,
which shows (CF3). Lastly, let
$f^\bullet,g^\bullet:X^\bullet\to Y^\bullet$ be any two morphisms
in $\Hot^?(\cA)$, and $s^\bullet:Y^\bullet\to Z^\bullet$ a
quasi-isomorphism such that
$s^\bullet\circ f^\bullet=s^\bullet\circ g^\bullet$ in
$\Hot^?(\cA)$. We rewrite the latter condition as
$s^\bullet\circ(f^\bullet-g^\bullet)=0$, and we notice the
exact sequence
$$
H^{-1}\rHot_\cA^\bullet(X^\bullet,(\Cone\,s)^\bullet)\to
H^0\rHot_\cA^\bullet(X^\bullet,Y^\bullet)\to
H^0\rHot^\bullet_\cA(X^\bullet,Z^\bullet)
$$
induced, via remark \ref{rem_long-exact-triang}(i), by the
distinguished triangle $\rHot^\bullet_\cA(X^\bullet,\Theta(s^\bullet))$
provided by proposition \ref{prop_triangul-hot}. Taking
into account example \ref{ex_Hom-complex-hots}(i), we
deduce that there exists a morphism
$$
h^\bullet:X^\bullet\to(\Cone\,s^\bullet)[-1]
\quad\text{such that}\quad
\partial_s^\bullet[-1]\circ h^\bullet=f^\bullet-g^\bullet
\qquad
\text{in $\Hot^?(\cA)$}
$$
where $\partial_s^\bullet:(\Cone\,s)^\bullet\to Y^\bullet[1]$ is
the boundary morphism of the true triangle $\Theta(s^\bullet)$.
Let now $\partial^\bullet_h:(\Cone\,h)^\bullet\to X^\bullet[1]$ be
the boundary morphism of $\Theta(h^\bullet)$; from proposition
\ref{prop_rotate-triangles} it follows easily that
$$
h^\bullet\circ\partial^\bullet_h[-1]=0
\qquad
\text{in $\Hot^?(\cA)$}
$$
whence $(f^\bullet-g^\bullet)\circ\partial^\bullet_h[-1]=0$,
and to deduce (CF4), it suffices to remark :

\begin{claim} $\partial^\bullet_h$ is a quasi-isomorphism.
\end{claim}
\begin{pfclaim}[] By assumption, $s^\bullet$ is a
quasi-isomorphism; in light of the long exact sequence
of remark \ref{rem_long-exact-triang}(i), it follows
easily that $H^i(\Cone\,s)^\bullet=0$ for every $i\in\Z$.
By the same token, we see that $H^i\partial^\bullet_h$
is an isomorphism for every $i\in\Z$, which is the claim.
\end{pfclaim}
\end{proof}

\sset\subsubsection{}\label{subsec_enough-objects}
Let now $\cA$ be an abelian category, $\cB$ an additive
subcategory of $\cA$, and set
$$
\Sigma_{\cB,+}:=\Sigma\cap\rMorph(\Hot^+(\cB))
\qquad
\Sigma_{\cB,-}:=\Sigma\cap\rMorph(\Hot^-(\cB))
$$
where $\Sigma\subset\rMorph(\Hot(\cA))$ is the set of all
quasi-isomorphisms. We may define the categories
$$
\sD^+_\cB(\cA):=\Hot^+(\cB)[\Sigma^{-1}_{\cB,+}]
\qquad
\sD^-_\cB(\cA):=\Hot^-(\cB)[\Sigma^{-1}_{\cB,-}]
$$
and clearly the inclusion functor $\cB\to\cA$ induces
natural functors
$$
i^+:\sD^+_\cB(\cA)\to\sD^+(\cA)
\qquad
i^-:\sD^-_\cB(\cA)\to\sD^-(\cA).
$$

\begin{proposition}\label{prop_enough-objects}
In the situation of \eqref{subsec_enough-objects},
suppose additionally that the following holds :
\begin{enumerate}
\alphaenu
\item
$\cB$ is a full subcategory of $\cA$.
\item
For every $A\in\Ob(\cA)$, there exists $B\in\Ob(\cB)$
with a monomorphism $A\to B$ in $\cA$.
\end{enumerate}
Then the functor $i^+$ is an equivalence, and
$\Sigma_{\cB,+}$ admits a left calculus of fractions.
\end{proposition}
\begin{proof} Set $\Sigma_+:=\Sigma\cap\rMorph(\Hot^+(\cA))$.
By theorem \ref{th_from-hot-to-der}(i), it suffices to show
that the induced functor
$$
\Hot^+(\cB)[\Sigma^{-1}_{\cB,+}]\to
\Hot^+(\cA)[\Sigma^{-1}_+]
$$
is an equivalence. In light of remark \ref{rem_full-subcat-fract}
and theorem \ref{th_from-hot-to-der}(ii), we are then reduced to
checking :

\begin{claim}\label{cl_enough-objects}
Suppose that conditions (a) and (b) of the proposition hold.
Then, for every $K^\bullet\in\sC^+(\cA)$ there exists
$L^\bullet\in\sC^+(\cB)$ with a quasi-isomorphism
$K^\bullet\to L^\bullet$.
\end{claim}
\begin{pfclaim}[] Say that $K^\bullet\in\sC^{\geq a}(\cA)$
for some $a\in\Z$; we shall construct first a double
complex $A^{\bullet\bullet}\in\sC^{\geq a}(\sC^{\geq-1}(\cA))$
such that
\begin{itemize}
\item
$A^{p,-1}=K^p$ and $d^{p-1,-1}_h=d_K^{p-1}$ for every $p\in\Z$.
\item
$A^{pq}\in\Ob(\cB)$ for every $p\in\Z$ and every $q\geq 0$.
\item
For every $p\in\Z$, the cohomology of the complex
$(A^{p\bullet},d_v^{p\bullet})$ vanishes in every degree.
\end{itemize}
To this aim, we proceed by induction on $p$. Hence,
we let $A^{pq}=0$ for every $p<a$ and every $q\in\Z$.
Suppose that $p\geq a$, and that $A^{ij}$ is already
defined for every $i<a$ and every $j\in\Z$, as well
as all the differentials $d^{ij}_v$ and $d^{i-1,j}_h$
for every $i<a$ and every $j\in\Z$. Then we construct
$A^{pj}$ and the differentials $d^{p,j-1}_v$, $d^{p-1,j}_h$
by induction on $j\in\Z$. Namely, for $j<-1$ we set
$A^{pj}=0$, and for $j=-1$ we set $A^{p,-1}=K^p$ and
$d^{p-1,-1}_h=d_K^{p-1}$. Thus, suppose that $j\geq 0$, and
$A^{p,j-1}$ is already given, as well as $d^{p-1,j-1}_h$
and $d^{p,j-2}_h$; we define the object $B^{pj}$ of $\cA$
as the coproduct in the push-out diagram
$$
\xymatrix{ A^{p-1,j} \ar[rr]^-g & & B^{pj} \\
\Coker\,d^{p-1,j-2}_v \ar[u]^{\bar d_v}
\ar[rr]^-{\bar d_h} & & \Coker\,d^{p,j-2}_v \ar[u]_f 
}$$
where $\bar d_v$ and $\bar d_h$ are induced by
$d^{p-1,j-1}_v$ and respectively $ d^{p-1,j-1}_h$.
Notice that, by inductive assumption, $\bar d_v$
is a monomorphism, so the same holds for $f$ (details
left to the reader). By assumption (b), we may find
a monomorphism $h:B^{pq}\to C$ for some $C\in\Ob(\cB)$,
and we set $A^{pq}:=C$, $d^{p-1,j}_h:=h\circ g$ and
$d_v^{p,j-1}:=h\circ f\circ p$, where
$p:A^{p,j-1}\to\Coker\,d^{p,j-2}_v$ is the natural
projection. Then it is easily seen that
$d^{p-1,j}_h\circ d^{p-1,j-1}_v=d^{p,j-1}_v\circ d^{p-1,j-1}_h$,
and moreover
$$
\Ker\,d^{p,j-1}_v/\Img\,d^{p,j-2}_v=0
$$
as required. Now, let
$t^{\geq 0}:\sC_2(\cA)\to\sC(\sC^{\geq 0}(\cA))$ be the
brutal truncation functor, and set
$$
B^{\bullet\bullet}:=t^{\geq 0}A^{\bullet\bullet}
\qquad
C^{\bullet\bullet}:=K^\bullet[0]
\qquad
L^\bullet:=\Tot\,B^{\bullet\bullet}.
$$
(Hence, $C^{\bullet\bullet}$ is the double complex whose rows
$C^{\bullet q}$ are the zero complex, except for $q=0$, where
it agrees with $K^\bullet$).
The differentials $d_v^{\bullet,-1}$ of $A^{\bullet\bullet}$
induce a morphism $f^\bullet:K^\bullet\to L^\bullet$, and it
remains only to check that $f^\bullet$ is a quasi-isomorphism.
To this aim, notice that the differentials $d_v^{\bullet,-1}$
also induce a morphism of double complexes
$$
g^{\bullet\bullet}:C^{\bullet\bullet}\to B^{\bullet\bullet}
$$
such that $\Tot\,g^{\bullet\bullet}=f^\bullet$. Moreover,
we have standard convergent spectral sequences
$$
E_2^{pq}:=H^q(C^{p,\bullet})\Rightarrow H^{p+q}\Tot\,C^{\bullet\bullet}
\qquad
F_2^{pq}:=H^q(B^{p,\bullet})\Rightarrow H^{p+q}\Tot\,B^{\bullet\bullet}
$$
and by construction, $g^{\bullet\bullet}$ induces an
isomorphism of spectral sequences
$E^{\bullet\bullet}_2\isom F^{\bullet\bullet}_2$, whence the claim.
\end{pfclaim}
\end{proof}

\begin{remark}\label{rem_enough-is-enough}
(i)\ \ 
In view of remark \ref{rem_derived-cat}(iii), we see
that the dual of proposition \ref{prop_enough-objects}
also holds. Namely, in the situation of
\eqref{subsec_enough-objects}, suppose that
\begin{enumerate}
\alphaenu
\item
$\cB$ is a full subcategory of $\cA$.
\item
For every $A\in\Ob(\cA)$, there exists $B\in\Ob(\cB)$
with an epimorphism $B\to A$ in $\cA$.
\end{enumerate}
Then the functor $i^-$ is an equivalence, and
$\Sigma_{\cB,-}$ admits a right calculus of fractions.

(ii)\ \
Two special cases of proposition \ref{prop_enough-objects}
(and its dual) are especially important. Namely, let
$\cI$ (resp. $\cP$) be the full subcategory of $\cA$
whose objects are the injective (resp. projective)
objects of $\cA$. We say that $\cA$ has {\em enough
injectives} (resp. {\em enough projectives}) if $\cI$
(resp. $\cP$) satisfies condition (b) of proposition
\ref{prop_enough-objects} (resp. condition (b) of (i)).
When this holds, proposition \ref{prop_enough-objects}
(and its dual) can be further sharpened; namely, we have :
\end{remark}

\begin{theorem}\label{th_enough-is-enough}
Let $\cA$ be any abelian category, and define the subcategories
$\cI$ and $\cP$ of $\cA$ as in remark
{\em\ref{rem_enough-is-enough}(ii)}. Then the following holds :
\begin{enumerate}
\item
Let $I^\bullet$ be any object of $\sC^+(\cI)$.
Let also $L^\bullet$ be any object of $\sC(\cA)$, and
$\alpha^\bullet:I^\bullet\to L^\bullet$ any quasi-isomorphism.
Then there exists a morphism $\beta:L^\bullet\to I^\bullet$
in $\sC(\cA)$ such that
$$
\beta^\bullet\circ\alpha^\bullet=\one_{I^\bullet}
\qquad
\text{in $\Hot(\cA)$}.
$$
\item
If $\cA$ has enough injectives, the inclusion functor
$\cI\to\cA$ induces an equivalence
$$
\Hot^+(\cI)\isom\sD^+(\cA).
$$
\item
Dually, if $\cA$ has enough projectives, the inclusion
$\cP\to\cA$ induces an equivalence
$$
\Hot^-(\cP)\isom\sD^-(\cA).
$$
\end{enumerate}
\end{theorem}
\begin{proof}(i): To begin with, we remark :

\begin{claim}\label{cl_null-homotopic}
In the situation of (i), let also $K^\bullet$ be any
acyclic complex of $\cA$, and
$\phi^\bullet:K^\bullet\to I^\bullet$ any morphism in
$\sC(\cA)$. Then $\phi^\bullet$ is null-homotopic. 
\end{claim}
\begin{pfclaim} We need to exhibit a system of morphisms
$(s^n:K^n\to I^{n-1}~|~n\in\Z)$ with
\set\begin{equation}\label{eq_null-homotopy}
\phi^n=d^{n-1}_I\circ s^n+s^{n+1}\circ d^n_K
\end{equation}
for every $n\in\Z$.
To this aim, say that $I^\bullet\in\Ob(\sC^{\geq a}(\cA))$
for some $a\in\Z$; then we let $s^j$ be the zero morphism
for every $j\leq a$. Next, suppose that $j\geq a$, and that
$s^j$ has already been given, so that \eqref{eq_null-homotopy}
holds with $n=j-1$. We compute :
$$
\begin{aligned}
(\phi^j-d^{j-1}_I\circ s^j)\circ d^{j-1}_K
=\, & d^{j-1}_I\circ\phi^{j-1}-d^{j-1}_I\circ s^j\circ d^{j-1}_K \\
=\, & d^{j-1}_I\circ(\phi^{j-1}-s^j\circ d^{j-1}_K) \\
=\, & d^{j-1}_I\circ d^{j-2}_I\circ s^{j-1} \\
=\, & 0.
\end{aligned}
$$
Since $\Img\,d^{j-1}_K=\Ker\,d^j_K$, it follows that
$\phi^j-d^{j-1}_I\circ s^j$ factors through a morphism
$\Img\,d^j_K\to I^j$, and since $I^j$ is injective,
the latter can be extended to a morphism
$s^{j+1}:K^{j+1}\to I^j$. By construction,
\eqref{eq_null-homotopy} holds for $n=j$, with this
choice of $s^{j+1}$, whence the claim.
\end{pfclaim}

Now, under the condition of (i), the complex
$\Cone\,\alpha^\bullet$ is acyclic (remark
\ref{rem_long-exact-triang}(iii)), so the boundary
morphism
$\partial^\bullet:\Cone\,\alpha^\bullet\to I^\bullet[1]$
is null-homotopic, by claim \ref{cl_null-homotopic}.
Let us then fix a homotopy $s^\bullet$ from the zero
morphism to $\partial^\bullet$; in each degree $j\in\Z$,
the morphism $s^j$ is the sum of a morphism
$s^j_L:L^j\to I^j$ and a morphism $s^j_I:I^{j+1}\to I^j$,
and the relation
$$
-\partial^j=d^{j-1}_{I[1]}\circ s^j+s^{j+1}\circ d^j_{\Cone\,\phi}
\qquad
\text{for every $j\in\Z$}
$$
translates as the pair of identities :
$$
\begin{aligned}
-d^j_I\circ s^j_L+s^{j+1}_L\circ d^j_L=\, & 0 \\
-d_I^j\circ s_I^j+s^{j+1}_L\circ\alpha^{i+1}+s^{j+1}_I\circ d^{i+1}_I
=\, & \one_{I^{j+1}}
\end{aligned}
$$
the first of which says that the system $(s^j_L~|~j\in\Z)$
amounts to a morphism of complexes
$s^\bullet_L:L^\bullet\to I^\bullet$, and the second says that
the system $(s^j_I~|~j\in\Z)$ yields a homotopy from
$\one_I$ to $s^\bullet_L\circ\alpha^\bullet$. In other
words, the morphism $\beta^\bullet:=s^\bullet_L$ will do.

(ii): By proposition \ref{prop_enough-objects}, it suffices
to check that every quasi-isomorphism in $\Hot^+(\cI)$ is
an isomorphism. But this follows immediately from (i).

(iii): Clearly the subcategory of injective objects
of $\cA^o$ is $\cP^o$, so the assertion follows from (ii)
and remark \ref{rem_derived-cat}(iii).
\end{proof}

\begin{corollary}\label{cor_enough-is-enough}
With the notation of theorem {\em\ref{th_enough-is-enough}},
the natural map
$$
\rHot_\cA(K^\bullet,I^\bullet)\to\Hom_{\sD(\cA)}(K^\bullet,I^\bullet)
$$
is an isomorphism, for every $K^\bullet\in\Ob(\sC(\cA))$ and
every $I^\bullet\in\Ob(\sC^+(\cI))$.
\end{corollary}
\begin{proof} By virtue of theorem \ref{th_from-hot-to-der},
every morphism $g^\bullet:K^\bullet\to I^\bullet$ in $\sD(\cA)$
is represented by a fraction $(f^\bullet,\alpha^\bullet)$,
where $f^\bullet:K^\bullet\to L^\bullet$ is any morphism
in $\Hot(\cA)$, and $\alpha^\bullet:I^\bullet\to L^\bullet$
is a quasi-isomorphism. Then, by theorem \ref{th_enough-is-enough}(i)
we may find a morphism $\beta^\bullet:L^\bullet\to I^\bullet$
in $\Hot(\cA)$ such that
$\beta^\bullet\circ\alpha^\bullet=\one_{I^\bullet}$; clearly
$\beta^\bullet$ is also a quasi-isomorphism, and
$g^\bullet=\beta^\bullet\circ f^\bullet$ in $\sD(\cA)$,
so the map of the corollary is surjective. For the
injectivity, suppose that
$f^\bullet_1,f^\bullet_2:K^\bullet\to I^\bullet$ are any
two morphisms in $\Hot(\cA)$ whose images agree in
$\sD(\cA)$; invoking again theorem \ref{th_from-hot-to-der},
we then find a quasi-isomorphism
$\alpha^\bullet:I^\bullet\to L^\bullet$ in $\Hot(\cA)$ such that
$\alpha^\bullet\circ f_1^\bullet=\alpha^\bullet\circ f_2^\bullet$
in $\Hot(\cA)$. Then pick again $\beta^\bullet$ as in the
foregoing; we deduce that
$$
f_1^\bullet=\beta^\bullet\circ\alpha^\bullet\circ f^\bullet_1=
\beta^\bullet\circ\alpha^\bullet\circ f^\bullet_2=f^\bullet_2
\qquad
\text{in $\Hot(\cA)$}
$$
whence the assertion.
\end{proof}

\begin{definition}\label{def_der-functors}
Let $\cA$ and $\cB$ be two abelian categories, and
$F:\sC^+(\cA)\to\sD(\cB)$ a functor. Denote
$\omega^+_\cA:\sC^+(\cA)\to\sD^+(\cA)$ and
$\omega^-_\cA:\sC^-(\cA)\to\sD^-(\cA)$ the localization
functors.
\begin{enumerate}
\item
A {\em right derived functor} of $F$ is a pair $(RF,\mu)$
consisting of a functor
$$
RF:\sD^+(\cA)\to\sD(\cB)
$$
and a natural transformation
$$
\mu:F\Rightarrow RF\circ\omega^+_\cA
$$
which satisfies the following universal property. For
every other pair $(G,\zeta)$ consisting of a functor
$G:\sD^+(\cA)\to\sD(\cB)$ and a natural transformation
$\zeta:F\Rightarrow G\circ\omega^+_\cA$, there exists
a unique natural transformation $\xi:RF\Rightarrow G$
such that
\set\begin{equation}\label{eq_univ-prop-RF}
\zeta=(\xi*\omega^+_\cA)\circ\mu.
\end{equation}
\item
Dually, if $F:\sC^-(\cA)\to\sD(\cB)$ is any functor, a
{\em left derived functor} of $F$ is a pair $(LF,\mu)$
consisting of a functor
$$
LF:\sD^-(\cA)\to\sD(\cB)
$$
and a natural transformation
$$
\mu:F\Rightarrow LF\circ\omega^-_\cA
$$
such that $((LF)^o,\mu^o)$ is a right derived functor
of $F^o$, under the natural identifications of remark
\ref{rem_derived-cat}(iii).
\item
Especially, let $\phi:\cA\to\cB$ be any additive functor,
and for ``$?$'' equal to $+$ or $-$, let
$\omega^?_\cB:\sC^?(\cB)\to\sD^?(\cB)$ be the localization
functor. Then a right (resp. left) derived functor of
$\omega^+_\cB\circ\sC^+(\phi)$ (resp. of
$\omega^-_\cB\circ\sC^-(\phi)$) shall just be called a
{\em right} (resp. {\em left}) {\em derived functor} of
$\phi$, and we shall just write $R\phi$ (resp. $L\phi$)
to denote this functor.
\end{enumerate}
\end{definition}

\begin{remark}\label{rem_acyclic-crit}
Keep the situation of definition \ref{def_der-functors}.

(i)\ \
As usual, if the right derived functor of $F$ exists, it is
unique up to isomorphism. More precisely, let $(RF,\mu)$ be
any pair as in definition \ref{def_der-functors}(i); if
$(G,\zeta)$ is any other such pair fulfilling the same
universal condition, there exists a unique isomorphism
of functors $\xi:RF\isom G$ such that
\eqref{eq_univ-prop-RF} holds. Likewise one characterizes
left functors up to isomorphism.

(ii)\ \
For $F$ as in definition \ref{def_der-functors}(i) (resp.
definition \ref{def_der-functors}(ii)), we shall use, for
every $p\in\Z$ and every object $K^\bullet$ of $\sD^+(\cA)$
(resp. of $\sD^-(\cA)$), the standard notation :
$$
R^pFK^\bullet:=H^p(RFK^\bullet)
\qquad
\text{(resp. $L^pFK^\bullet:=H^p(LFK^\bullet)$)}.
$$
Also, if $\phi$ is as in definition \ref{def_der-functors}(iii),
and the right (resp. left) derived functor of $\phi$ exists,
then we have a natural transformation
$$
\phi A\to R^0\phi A[0]
\qquad
\text{(resp.\ $L^0\phi A[0]\to\phi A$)}
\qquad
\text{for every $A\in\Ob(\cA)$}
$$
and it is easily seen that this transformation is an
isomorphism of functors if and only $\phi$ is left exact
(resp. right exact : details left to the reader).

(iii)\ \
Suppose that $\cA$ has enough injectives (see remark
\ref{rem_enough-is-enough}(ii)), and $F:\sC^+(\cA)\to\sD(\cB)$
is a functor that factors through $\Hot^+(\cA)$ (via
the natural functor \eqref{eq_cplx-to-hot}). Then the
right derived functor of $F$ exists, and can be constructed
as follows. First, say that $f^\bullet:I^\bullet\to J^\bullet$
is a quasi-isomorphism, with
$I^\bullet,J^\bullet\in\Ob(\sC^+(\cI))$ (where $\cI$ is
as in remark \ref{rem_enough-is-enough}(ii)). By theorem
\ref{th_enough-is-enough}(i), it follows that $f^\bullet$
is a homotopy equivalence, so $Ff^\bullet$ is an isomorphism
by assumption, and therefore $F$ induces a functor
$$
RF_\cI:\sD^+_\cI(\cA)\to\sD^+(\cB)
$$
(notation of \eqref{subsec_enough-objects}). However,
according to proposition \ref{prop_enough-objects}, the
natural functor $i^+:\sD^+_\cI(\cA)\to\sD^+(\cA)$ admits a
quasi-inverse $j^+:\sD^+(\cA)\to\sD^+_\cI(\cA)$; in view
of claim \ref{cl_enough-objects}, the functor $j^+$ can
be described by choosing, for every object $K^\bullet$
of $\sC^+(\cA)$, a quasi-isomorphism
$$
\tau_K^\bullet:K^\bullet\to I_K^\bullet
\qquad
\text{with $I_K^\bullet\in\Ob(\sC^+(\cI))$}.
$$
With this notation, we let
$$
RF:=RF_\cI\circ j^+
$$
and we define a natural transformation
$\mu:F\Rightarrow RF\circ\omega^+_\cA$, by the rule :
$$
K^\bullet\mapsto F(\tau_K^\bullet)
\qquad
\text{for every $K^\bullet\in\Ob(\sC^+(\cA))$}.
$$
Now, suppose that $G:\sD^+(\cA)\to\sD(\cB)$ is another
functor, with a natural transformation
$\zeta:F\Rightarrow G\circ\omega^+_\cA$. There follows
a commutative diagram in $\sD(\cB)$
$$
\xymatrix{ FK^\bullet \ar[r]^-{\zeta^\bullet_K}
\ar[d]_{F\tau^\bullet_K} & GK^\bullet \ar[d]^{G\tau^\bullet_K} \\
FI^\bullet_K \ar[r]^-{\zeta^\bullet_{I_K}} & GI^\bullet_K
}$$
in which $G\tau^\bullet_K$ is an isomorphism. It follows
that \eqref{eq_univ-prop-RF} holds if and only if
$\xi_K=(G\tau^\bullet_K)^{-1}\circ\zeta^\bullet_{I_K}$ for
every $K^\bullet\in\Ob(\sC^+(\cA))$, and therefore $(RF,\mu)$
is a derived functor of $F$, as stated. This is essentially
the original construction of the right derived functor
proposed by Grothendieck in \cite{Toh}, which predates
the invention of derived categories. Dually, one
obtains likewise a left derived functor of $F$, in case
$\cA$ has enough projectives : in this case, one starts
by fixing, for every bounded above complex $K^\bullet$
of $\cA$, a quasi-isomorphism $P^\bullet_K\to K^\bullet$
with $P^\bullet$ in $\sC^+(\cP)$.

(iv)\ \
From \cite[Lemme 3.3.1]{Toh} one can also extract a
construction which works only for additive functors
$F:\cA\to\cB$, but which requires only that $\cA$ has
{\em enough $F$-acyclic objects}. Namely, let us suppose
that $\cA$ contains a full subcategory $\cM$ with the
following properties :
\begin{enumerate}
\alphaenu
\item
For every $A\in\mathrm{Ob}(\cA)$, there exists a monomorphism
$A\to M$, with $M\in\mathrm{Ob}(\cM)$.
\item
For every $M,M'\in\Ob(\cM)$ and every short exact sequence in $\cA$
$$
\Sigma\qquad :\qquad
0\to M'\to M\to M''\to 0
\qquad\qquad\qquad
$$
the object $M''$ is also in $\cM$, and $\Sigma$ induces
a short exact sequence
$$
F(\Sigma)\qquad :\qquad
0\to F(M')\to F(M)\to F(M'')\to 0.
\qquad\qquad\qquad
$$
\romanenu
\end{enumerate}
Then we obtain a right derived functor of $F$ as
follows. First, let $(M^\bullet,d^\bullet_M)$ be any
acyclic object of $\sC^+(\cM)$; an easy induction
shows that the induced sequences
$$
\Sigma^i
\quad :\quad
0\to\Img\,d^i\to M^{i+1}\to\Img\,d^{i+1}\to 0
$$
are short exact and $\Img\,d^i\in\Ob(\cM)$ for every
$i\in\Z$. By assumption, the sequences $F(\Sigma^i)$
are then also short exact, and therefore $FM^\bullet$
is still acyclic. Next, let $f^\bullet:M^\bullet\to N^\bullet$
be any quasi-isomorphism in $\sC^+(\cM)$; then
$\Cone\,f^\bullet$ is acyclic, so the same holds for
$F(\Cone\,f^\bullet)$, by the foregoing, and since
$\sC^+(F)$ is a triangulated functor, we conclude
that $Ff^\bullet$ is still a quasi-isomorphism. Thus,
$\sC^+(F)$ induces a functor
$$
RF_\cM:\sD^+_\cM(\cA)\to\sD^+(\cB).
$$
On the other hand, proposition \ref{prop_enough-objects}
says that the natural functor $i^+:\sD_\cM^+(\cA)\to\sD^+(\cA)$
admits a quasi-inverse $j^+:\sD^+(\cA)\to\sD^+_\cM(\cA)$,
so we may set again
$$
RF:=RF_\cM\circ j^+
$$
and one obtains as in (iii) a natural transformation
$\mu$, such that the pair $(RF,\mu)$ fulfills the required
universal property : details left to the reader.

(v)\ \
Notice that any short exact sequence
$0\to I'\to I\to I''\to 0$ in $\cA$ with $I'$ injective,
is split, and $I''$ is injective if and only if the same
holds for $I$. Hence, if $\cA$ has enough injectives, the
category $\cM:=\cI$ fulfills the conditions of (iv) for
every additive functor $F$. Conversely, suppose that
$\cM$ fulfills the conditions of (iv), and moreover,
every $A\in\mathrm{Ob}(\cC)$ isomorphic to a direct
factor of an object of $\cM$, is an object of $\cM$.
Then $\cM$ contains every injective object of $\cA$.
Indeed, if $I$ is an injective object of $\cC$, then
by (a) we can find a monomorphism $I\to M$ with $M$
in $\cM$; hence $I$ is a direct summand of $M$, so it
is in $\cM$, by our assumption.

(vi)\ \
Notice that the derived functor $RF$ constructed in (iv)
is triangulated (see remark \ref{rem_long-exact-triang}(ii));
moreover, the construction implies immediately that
$$
R^pFM[0]=0
\qquad
\text{for every $p\neq 0$}
$$
whenever $M$ is an $F$-acyclic object of $\cA$.
\end{remark}

\begin{example}\label{ex_construct-RHom}
Let $\cA$ be any abelian category with enough injectives;
recall that the functor
$$
\Hom^\bullet_\cA:\sC(\cA)\times\sC(\cA)^o\to\sC(\Z\Mod)
$$
factors through $\Hot(\cA)\times\Hot(\cA)^o$ (see example
\ref{ex_hot-morphisms}); hence, after fixing an
equivalence $j^+:\sD^+(\cA)\to\sD^+_\cI(\cA)$ and arguing
as in remark \ref{rem_acyclic-crit}(iii), we deduce a
functor
$$
R\Hom_\cA^\bullet:\sD^+(\cA)\times\Hot(\cA)^o\to\sD(\Z\Mod)
$$
such that, for every $K^\bullet\in\Ob(\Hot(\cA))$, the
restriction
\set\begin{equation}\label{eq_restrict-RHom}
\sD^+(\cA)\to\sD(\Z\Mod)
\quad :\quad
L^\bullet\mapsto R\Hom^\bullet_\cA(K^\bullet,L^\bullet)
\end{equation}
is the derived functor of the functor
$\Hom^\bullet_\cA(K^\bullet,-):\sC^+(\cA)\to\sD(\Z\Mod)$.
Moreover, \eqref{eq_restrict-RHom} is triangulated, by
proposition \ref{prop_triangul-hot}. Furthermore, let
$\phi^\bullet:K^\bullet_1\to K^\bullet_2$ be any quasi-isomorphism
in $\Hot(\cA)$; we claim that
$R\Hom^\bullet_\cA(\phi^\bullet,L^\bullet)$ is a quasi-isomorphism,
for every $L^\bullet\in\Ob(\sD^+(\cA))$. Indeed, set
$I^\bullet_L:=j^+L^\bullet$; by example \ref{ex_Hom-complex-hots}(i)
and corollary \ref{cor_enough-is-enough} we have natural isomorphisms
of abelian groups
$$
H^nR\Hom^\bullet_\cA(\phi^\bullet,L^\bullet)\isom
\rHot_\cA(\phi^\bullet,I^\bullet_L[n])\isom
\Hom_{\sD(\cA)}(\phi^\bullet,I^\bullet_L[n])
$$
whence the contention. It follows that $R\Hom^\bullet_\cA$
descends to an additive functor
$$
R\Hom^\bullet_\cA:\sD^+(\cA)\times\sD(\cA)^o\to\sD(\Z\Mod)
$$
with a natural isomorphism of abelian groups, for every
$K^\bullet$ in $\sD(\cA)$ and every $L^\bullet$ in $\sD^+(\cA)$ :
\set\begin{equation}\label{eq_shift-and-shout}
R^n\Hom^\bullet_\cA(K^\bullet,L^\bullet)\isom
\Hom_{\sD(\cA)}(K^\bullet,L^\bullet[n]).
\end{equation}
\end{example}

\begin{remark}\label{rem_long-exact-derHom}
In the situation of example \ref{ex_construct-RHom}, our
methods do not allows us to extend the functor
$R\Hom^\bullet_\cA$ to unbounded complexes in both arguments,
but still we can show that every $L^\bullet\in\sD(\cA)$ and
every distinguished triangle of $\sD(\cA)$
$$
\Theta
\quad :\quad
K^\bullet\to K'^\bullet\to K''^\bullet\to K^\bullet[1]
$$
induce a long exact $\Ext$-sequence :
$$
\Hom_{\sD(\cA)}(L^\bullet,K^\bullet)\!\to\!
\Hom_{\sD(\cA)}(L^\bullet,K'^\bullet)\!\to\!
\Hom_{\sD(\cA)}(L^\bullet,K''^\bullet)\!\to\!
\Hom_{\sD(\cA)}(L^\bullet,K^\bullet[1])\!\to\!\cdots
$$
Indeed, let $L'^\bullet\to L^\bullet$ be any quasi-isomorphism;
we may assume that $\Theta$ is already a distinguished triangle
of $\Hot(\cA)$, and due to proposition \ref{prop_triangul-hot}
and remark \ref{rem_long-exact-triang}(iii), we get the long
exact sequence
$$
\rHot_\cA(L'^\bullet,K^\bullet)\to
\rHot_\cA(L'^\bullet,K'^\bullet)\to
\rHot_\cA(L'^\bullet,K''^\bullet)\to
\rHot_\cA(L'^\bullet,K^\bullet[1])\to\cdots
$$
and on the other hand, theorem \ref{th_from-hot-to-der}(ii)
implies that $\Hom_{\sD(\cA)}(L^\bullet,K^\bullet)$ is isomorphic
to the colimit of the system
$(\rHot_\cA(L'^\bullet,K^\bullet)~|~L'^\bullet\to L^\bullet)$
indexed by the filtered set of all such quasi-isomorphisms
(proposition \ref{prop_calculus-frac}(i)). Since all filtered
colimits are exact in the category of abelian groups, the
assertion follows. Likewise, from $\Theta$ and any
$L^\bullet\in\Ob(\sD(\cA))$ we get as well the long exact
$\Ext$-sequence :
$$
\Hom_{\sD(\cA)}(K^\bullet[1],L^\bullet)\!\to\!
\Hom_{\sD(\cA)}(K''^\bullet,L^\bullet)\!\to\!
\Hom_{\sD(\cA)}(K'^\bullet,L^\bullet)\!\to\!
\Hom_{\sD(\cA)}(K^\bullet,L^\bullet)\!\to\!\cdots
$$
\end{remark}

\begin{example}\label{ex_derive-biadditive}
(i)\ \
Remark \ref{rem_acyclic-crit}(iv) can be adapted to
biadditive functors. Namely, in the situation of
\eqref{subsec_hot-biadditive}, let us say that an
object $A$ of $\cA$ (resp. $A'$ of $\cA'$) is
{\em left $B$-flat} (resp. {\em right $B$-flat}) if the
additive functor $B(A,-)$ (resp. $B(-,A')$) is exact.
We assume that :
\begin{itemize}
\item
$\cA$ has {\em enough left $B$-flat objects} and $\cA'$
has {\em enough right $B$-flat objects}; {\em i.e.}
for every $M\in\Ob(\cA)$ (resp. $M'\in\Ob(\cA')$) there
exists an epimorphism $A\to M$ (resp. $A'\to M'$), where
$A$ (resp. $A'$) is a left (resp. right) $B$-flat object
of $\cA$ (resp. of $\cA'$).
\end{itemize}
Let us denote by $\cF$ (resp. $\cF'$) the full subcategory
of $\cA$ (resp. of $\cA'$) whose objects are the left
$B$-flat (resp. right $B$-flat) objects of $\cA$ (of $\cA'$).
By remark \ref{rem_enough-is-enough}(i), the natural functors
\set\begin{equation}\label{eq_flat-derived-obj}
\sD^-_\cF(\cA)\to\sD^-(\cA)
\qquad
\sD^-_{\cF'}(\cA')\to\sD^-(\cA')
\end{equation}
are equivalences. On the other hand, let $P^\bullet$ be any
object of $\sC^-(\cF)$, and $f^\bullet:K^\bullet\to L^\bullet$
any quasi-isomorphism in $\sC^-(\cA')$; we wish to show that
the induced morphism $B^\bullet_-(P^\bullet,f^\bullet)$ is a
quasi-isomorphism as well. By remark \ref{rem_biadditive-triang},
it suffices to check that $B^\bullet_-(P^\bullet,\Cone\,f^\bullet)$
is acyclic. However, the double complex
$B^{\bullet\bullet}(P^\bullet,\Cone\,f^\bullet)$ yields a
convergent spectral sequence
$$
E^{ij}_1:=H^iB(P^j,\Cone\,f^\bullet)\Rightarrow
H^{i+j}B^\bullet_-(P^\bullet,\Cone\,f^\bullet)
$$
and since $P^j$ is left $B$-flat and $\Cone\,f^\bullet$
is acyclic, we get $E^{ij}_1=0$ for every $i,j\in\Z$,
whence the claim. Likewise, we see that
$B^\bullet(g^\bullet,Q^\bullet)$ is a quasi-isomorphism,
whenever $Q^\bullet\in\Ob(\sC^-(\cF'))$ and $g^\bullet$
is any quasi-isomorphism in $\sC^-(\cA)$. Then, as in
remark \ref{rem_acyclic-crit}(iv), after fixing quasi-inverse
functors for the equivalences \eqref{eq_flat-derived-obj},
we obtain a functor
$$
LB^\bullet_-:\sD^-(\cA)\times\sD^-(\cA')\to\sD^-(\cA'')
$$
with a natural transformation
$$
\mu:LB^\bullet_-\circ(\omega^-_\cA\times\omega^-_{\cA'})
\to\omega^-_{\cA''}\circ B^\bullet_-
$$
such that, for every $K^\bullet\in\Ob(\sC^-(\cA))$,
the induced functor
$$
\sD^-(\cA')\to\sD^-(\cA'')
\quad :\quad
X^\bullet\mapsto LB^\bullet(K^\bullet,X^\bullet)
$$
together with the corresponding restriction of $\mu$,
is a left derived functor of $B^\bullet_-(K^\bullet,-)$.
Symmetrically, for every $K^\bullet\in\Ob(\sC^-(\cA'))$,
the functor $LB^\bullet_-(-,K^\bullet)$ provides a left
derived functor of $B^\bullet_-(-,K^\bullet)$ : details
left to the reader.

(ii)\ \
By remark \ref{rem_acyclic-crit}(ii), the following
conditions are equivalent :
\begin{enumerate}
\alphaenu
\item
The induced morphism
$$
B^*(A,A'):=L^0B^\bullet(A[0],A'[0])\to B(A,A')
$$
is an isomorphism for every $A\in\Ob(\cA)$ and every
$A'\in\Ob(\cA')$.
\item
The functor $B(A,-):\cA'\to\cA''$ is right exact for every
$A\in\Ob(\cA)$.
\item
The functor $B(-,A'):\cA\to\cA''$ is right exact for every
$A'\in\Ob(\cA)$.
\end{enumerate}
Moreover, $B^*$ is obviously another biadditive functor
$\cA\times\cA'\to\cA''$, and every right (resp. left)
$B$-flat object of $\cA$ (of $\cA'$) is also a right
(resp. left) $B^*$-flat object. Thus, $\cA$ has enough
right $B^*$-flat objects and $\cA'$ has enough left
$B^*$-flat objects, so we may define as well the left
defived functor $LB^*$ of $B^*$. Furthermore, a simple
inspection shows that the functors $B^*(A,-)$ and $B^*(-,A')$
are right exact, for every $A\in\Ob(\cA)$ and $A'\in\Ob(\cA')$,
and the induced morphism
$$
LB^*\to LB
$$
is an isomorphism of functors. Lastly, an object
$P$ of $\cA$ (resp. $P'$ of $\cA'$) is right (resp.
left) $B^*$-flat if and only if $L^1B(P,-)$ (resp.
$L^1B(-,P')$) is the (constant) zero functor. The
detailed verifications of all these assertions shall
be left as exercises for the reader.
\end{example}

\sset\subsubsection{}\label{subsec_give_Ps}
Let $(\cA,\otimes,\Phi,\Psi)$ be any abelian tensor category;
example \ref{ex_derive-biadditive} applies especially to
the functor $\otimes$. In this case, it is clear that
an object of $\cA$ is left $\otimes$-flat if and only
if it is right $\otimes$-flat, and such objects shall
therefore be called simply {\em $\otimes$-flat}. Thus,
suppose that $\cA$ has enough $\otimes$-flat objects;
we conclude that the induced tensor functor for complexes
of $\cA$ (see example \ref{ex_monoidal}) admits a left
derived functor
$$
-\derotimes-:\sD^-(\cA)\times\sD^-(\cA)\to\sD^-(\cA)
$$
called the {\em derived tensor product}. Explicitly, this
is defined as follows : for every bounded above complex
$K^\bullet$ we fix a quasi-isomorphism 
$\rho^\bullet_K:P_K^\bullet\to K^\bullet$ with $P_K^\bullet$
a bounded above complex of $\otimes$-flat objects, and
we set
$$
K^\bullet\derotimes L^\bullet:=P_K^\bullet\otimes P_L^\bullet.
$$
By inspecting the construction, we also see that
$K^\bullet\derotimes L^\bullet$ is naturally isomorphic --
via $\rho^\bullet_K$ and $\rho^\bullet_L$ -- to both
$K^\bullet\otimes P^\bullet_L$ and $P^\bullet_K\otimes L^\bullet$.
Moreover, if $A$ is any object of $\cA$, clearly the functor
$A\otimes-:\cA\to\cA$ is right exact if and only if the
same holds for the functor $-\otimes A:\cA\to\cA$, and
remark \ref{rem_acyclic-crit}(ii) implies that the latter
condition holds if and only if the induced morphism
$$
H^0(A[0]\derotimes B[0])\to A\otimes B
$$
is an isomorphism for every $B\in\Ob(\cA)$. We also let
$$
\Tor^\cA_i(K^\bullet,L^\bullet):=H_i(K^\bullet\derotimes L^\bullet)
\qquad
\text{for every $i\in\Z$}.
$$
In case $\cA=A\Mod$ for some ring $A$, it is customary to
denote this functor by $-\derotimes_A-$, and then one also
writes $\Tor^A_i$ instead of $\Tor^{A\Mod}_i$.

\begin{remark}\label{rem_der-tensor-varie}
(i)\ \
Using the commutativity and associativity constraints
for the tensor product in $\cA$, we deduce -- in light
of example \ref{ex_monoidal}(i,ii) -- natural
{\em associativity isomorphisms}
$$
K^\bullet\derotimes(L^\bullet\derotimes Q^\bullet)
\isom
(K^\bullet\derotimes L^\bullet)\derotimes Q^\bullet
\qquad
\text{in $\sD^-(\cA)$}
$$
as well as {\em commutativity isomorphisms}
$$
K^\bullet\derotimes L^\bullet\isom L^\bullet\derotimes K^\bullet
\qquad
\text{in $\sD^-(\cA)$}
$$
for any bounded above complexes $K^\bullet,L^\bullet$, and
$Q^\bullet$.

(ii)\ \
In the situation of \eqref{subsec_give_Ps}, take $\cA=A\Mod$
for some ring $A$, and suppose furthermore that $\phi:A\to B$
is a ring homomorphism, $K^\bullet$ a bounded above complex of
$A$-modules, and $L^\bullet$ a complex of $B$-modules. Notice
that $P^\bullet_K\otimes_AL^\bullet$ is naturally a complex of
$B$-modules; also, if $L^\bullet\to Q^\bullet$ is any
morphism of complexes of $B$-modules, then the induced
map
$P^\bullet_K\otimes_AL^\bullet\to P^\bullet_K\otimes_AQ^\bullet$
is $B$-linear. It follows easily that the derived tensor
product yields a functor
$$
\sD^-(A\Mod)\times\sD^-(B\Mod)\to\sD^-(B\Mod)
\qquad
(K^\bullet,L^\bullet)\mapsto K^\bullet\derotimes_AL^\bullet
$$
such that, denoting $\phi^*:\sD^-(B\Mod)\to\sD^-(A\Mod)$ the
``forgetful'' functor, we have a natural isomorphism
$$
K^\bullet\derotimes_A\phi^*L^\bullet\isom
\phi^*(K^\bullet\derotimes_AL^\bullet)
\qquad
\text{in $\sD^-(A\Mod)$}.
$$
\end{remark}

\sset\subsubsection{}\label{subsec_pairing-Tors}
Let now $M^\bullet_1$, $M^\bullet_2$, $N^\bullet_1$, $N^\bullet_2$
be any four objects of $\sC^-(\cA)$, and to ease notation, set
$$
M^\bullet_{12}:=M^\bullet_1\otimes M^\bullet_2
\qquad
N^\bullet_{12}:=N^\bullet_1\otimes N^\bullet_2
\qquad
P_{12}^\bullet:=P_{M_{12}}^\bullet
\qquad
\rho^\bullet_{12}:=\rho^\bullet_{M_{12}}
$$
as well as
$P^\bullet_i:=P^\bullet_{M_i}$ and $\rho_i^\bullet:=\rho^\bullet_{M_i}$
for $i=1,2$. There is a commutative diagram in $\sD^-(\cA)$
$$
\xymatrix{ P^\bullet_1\otimes P^\bullet_2 \ar[rr]^-{\phi^\bullet_{12}}
\ar[rd]_{\rho^\bullet_1\otimes\rho^\bullet_2} & &
P^\bullet_{12} \ar[ld]^{\rho^\bullet_{12}} \\
& M^\bullet_{12}
}$$
where $\rho^\bullet_{12}$ is an isomorphism, so $\phi^\bullet_{12}$
is uniquely determined, whence a map
$$
(M_1^\bullet\derotimes N_1^\bullet)\otimes(M^\bullet_2\derotimes N_2^\bullet)
\isom(P^\bullet_1\otimes P^\bullet_2)\otimes N^\bullet_{12}
\xrightarrow{\ \phi^\bullet_{12}\otimes N^\bullet_{12}\ }
P^\bullet_{12}\otimes N^\bullet_{12}\isom
M_{12}^\bullet\derotimes N^\bullet_{12}.
$$
Taking into account example \ref{ex_monoidal}(iii) and
remark \ref{rem_derived-cat}(v), we deduce a bilinear pairing
$$
\Tor^\cA_i(M^\bullet_1,N^\bullet_1)\otimes
\Tor^\cA_j(M^\bullet_2,N^\bullet_2)\to
\Tor^\cA_{i+j}(M^\bullet_{12},N^\bullet_{12})
\qquad
\text{for every $i,j\in\Z$}.
$$
Moreover, suppose that $M^\bullet_3$ and $N^\bullet_3$ are two
other bounded above complexes of $\cA$; by inspecting the
constructions, we find a commutative diagram
$$
\xymatrix{ P^\bullet_1\otimes(P^\bullet_2\otimes P^\bullet_3)
\ar[rrr]^-{P^\bullet_1\otimes\phi^\bullet_{23}}
\ar[rrrd]_{\rho^\bullet_1\otimes(\rho^\bullet_2\otimes\rho^\bullet_3)\ \ }
\ar[ddd]_{\Phi^\bullet_P} & & &
P^\bullet_1\otimes P^\bullet_{23}
\ar[rrr]^-{\phi^\bullet_{1,23}} \ar[d]^{\rho^\bullet_1\otimes\rho^\bullet_{23}}
& & &
P^\bullet_{1,23} \ar[llld]^{\rho^\bullet_{1,23}} \ar[ddd]^{P^\bullet_{\Phi_M}} \\
& & & M^\bullet_{1,23} \ar[d]^{\Phi^\bullet_M} \\
& & & M^\bullet_{12,3} \\
(P^\bullet_1\otimes P^\bullet_2)\otimes P^\bullet_3
\ar[rrr]^-{\phi^\bullet_{12}\otimes P^\bullet_3}
\ar[rrru]^{(\rho^\bullet_1\otimes\rho^\bullet_2)\otimes\rho^\bullet_3\ \ }
& & & P^\bullet_{12}\otimes P^\bullet_3
\ar[rrr]^-{\phi^\bullet_{12,3}} \ar[u]_{\rho^\bullet_{12}\otimes\rho^\bullet_3}
& & & P^\bullet_{12,3} \ar[lllu]_{\rho^\bullet_{12,3}}
}$$
where $M^\bullet_{1,23}:=M_{12}^\bullet\otimes M^\bullet_3$,
$M^\bullet_{23}:=M^\bullet_2\otimes M^\bullet_3$,
$M^\bullet_{1,23}:=M_1^\bullet\otimes M^\bullet_{23}$,
and likewise for $P^\bullet_{1,23}$, $P^\bullet_{23}$, and
$P^\bullet_{12,3}$ and the morphism $\phi^\bullet_{23}$,
$\phi^\bullet_{1,23}$, $\phi^\bullet_{12,3}$, $\rho^\bullet_{23}$,
$\rho^\bullet_{1,23}$, $\rho^\bullet_{12,3}$. Here $\Phi^\bullet_M$
and $\Phi^\bullet_P$ are the associativity constraints.

Therefore, set
$T^j_i:=\Tor^\cA_i(M^\bullet_j,N^\bullet_j)$ for every $i\in\Z$
and $j=1,2,3$, and also
$$
T_i^{jk}:=\Tor^\cA_i(M^\bullet_{jk},N^\bullet_{jk})
\qquad
T_i^{1,23}:=\Tor^\cA_i(M^\bullet_{1,23},N^\bullet_{1,23})
\qquad
T_i^{12,3}:=\Tor^\cA_i(M^\bullet_{12,3},N^\bullet_{12,3})
$$
for every $i\in\Z$, with $j=1,2$ and $k=j+1$; in light of example
\ref{ex_monoidal}(iv), we deduce a commutative diagram in $\cA$ :
\set\begin{equation}\label{eq_laborious}
{\diagram T_i^1\otimes(T^2_j\otimes T^3_k)
\ar[r] \ar[d] & T^1_i\otimes T_{j+k}^{23} \ar[r] &
T^{1,23}_{i+j+k} \ar[d] \\
(T^1_i\otimes T^2_j)\otimes T^3_k \ar[r] &
T_{i+j}^{12}\otimes T^3_k \ar[r] & T^{12,3}_{i+j+k}
\enddiagram}
\end{equation}
whose horizontal arrows are given by the above bilinear
pairing, and whose left (resp. right) vertical arrow is
the associativity constraint (resp. is induced
by the associativity constraint $\Phi^\bullet_M$).

The following lemma is borrowed from
\cite[Exp.VII,Prop.1.8(ii)]{SGA6} :

\begin{lemma}\label{lem_shift-and-shout}
Let $\cA$ be any abelian category and $a,b\in\Z$ any two
integers. We have :
\begin{enumerate}
\item
Let $K^\bullet\in\Ob(\sD^{\leq b}(\cA))$ and
$L^\bullet\in\Ob(\sD^{\geq a}(\cA))$ be any two complexes,
and suppose that $\cA$ has enough injectives. Then the
complex $R\Hom^\bullet_\cA(K^\bullet,L^\bullet)$ lies in
$\sD^{\geq a-b}(\cA)$.
\item
Suppose that $(\cA,\otimes\Phi,\Psi)$ is an abelian
tensor category with enough $\otimes$-flat objects,
and let $K^\bullet\in\Ob(\sD^{\leq a}(\cA))$ and
$L^\bullet\in\Ob(\sD^{\leq b}(\cA))$. then
$K^\bullet\derotimes L^\bullet\in\Ob(\sD^{\leq a+b}(\cA))$.
\end{enumerate}
\end{lemma}
\begin{proof}(i) follows easily from
\eqref{eq_shift-and-shout} and remark
\ref{rem_derived-cat}(iv) : details left to the reader.
Assertion (ii) follows immediately from the construction
of the derived tensor product.
\end{proof}

\sset\subsubsection{Complexes of modules over a ring}
We conclude this section with a discussion of a few special
features of the category of complexes of modules over a ring
and of its derived category.

\sset\subsubsection{}\label{subsec_ease-notation}
Let $A$ be any ring; to ease notation we set
$$
\sC(A):=\sC(A\Mod)
\qquad
\Hot(A):=\Hot(A\Mod)
\qquad
\sD(A):=\sD(A\Mod)
$$
and likewise we define $\sC^I(A)$, $\Hot^I(A)$ and $\sD^I(A)$
for any interval $I\subset\Z$.

\begin{proposition}\label{prop_depressed}
With the notation of \eqref{subsec_ease-notation}, the
following holds :
\begin{enumerate}
\item
The localization functor $j^I:\sC^I(A)\to\sD^I(A)$
commutes with direct sums and direct products.
\item
Let $(K_n^\bullet,\phi^\bullet_n:K^\bullet_n\to K^\bullet_{n+1}~|~n\in\N)$
be any direct system of objects of $\sC^I(A)$,  and
$L^\bullet$ its colimit in the category $\sC^I(A)$. Then
$L^\bullet$ represents also the colimit of the system
$(j^I(K^\bullet_n),j^I(\phi^\bullet_n)~|~n\in\N)$ in the
category $\sD^I(A)$.
\item
Let $(K_n^\bullet,\phi^\bullet_n:K^\bullet_{n+1}\to K^\bullet_n~|~n\in\N)$
be any inverse system of complexes of $A$-modules, and $L^\bullet$
its limit in the category $\sC(A)$. Suppose that for every
$i\in\Z$ the inverse system $(K_n^i,\phi^i_n~|~n\in\N)$ satisfies
the Mittag-Leffler condition. Then for every object $C^\bullet$ of
$\sD(A)$ we have a short exact sequence
$$
0\to\lim_{n\in\N}{}^{\!1}\,\Hom_{\sD(A)}(C^\bullet,K^\bullet_n[-1])\to
\Hom_{\sD(A)}(C^\bullet,L^\bullet)\to
\lim_{n\in\N}\Hom_{\sD(A)}(C^\bullet,K^\bullet_n)\to 0.
$$
\end{enumerate}
\end{proposition}
\begin{proof}(i): In light of remark \ref{rem_homotopies}(vi),
it suffices to check that the localization functor
$$
\Hot^I(A)\to\sD^I(A)
$$
commutes with direct sums and direct products. Hence, consider
any family $(K^\bullet_j~|~j\in J)$ of objects of $\Hot^I(A)$
indexed by a small set $J$, set $K^\bullet:=\prod_{j\in J}K_j^\bullet$,
the direct product in the category $\sC^I(A)$, and denote by
$\pi_j^\bullet:K^\bullet\to K^\bullet_j$ the canonical projection,
for every $j\in J$; we know that $K^\bullet$ represents also
the direct product of the family $K^\bullet_\bullet$ in
$\Hot^I(A)$. Hence, let $L^\bullet$ be any other object of
$\Hot^I(A)$; from any morphism
$\phi^\bullet:L^\bullet\to K^\bullet$ in $\sD^I(A)$ we obtain
the system $(\phi^\bullet_j:=\pi^\bullet_j\circ\phi^\bullet:
K^\bullet_j\to L^\bullet~|~j\in J)$ of morphisms in $\sD^I(A)$.
Conversely, given such a system of morphisms, we know by
theorem \ref{th_from-hot-to-der}(ii) that there exist for
every $j\in J$ a quasi-isomorphism
$\psi^\bullet_j:K'^\bullet_j\to K^\bullet_j$ and a morphism
$\tilde\phi^\bullet_j:K'^\bullet\to L^\bullet$ in $\sC^I(A)$
representing the class of $\phi^\bullet_j$ in $\Hot^I(A)$.
The systems $(\tilde\phi^\bullet_j~|~j\in J)$ and
$(\psi^\bullet_j~|~j\in J)$ yield morphisms
$\tilde\phi^\bullet:K'^\bullet:=\prod_{j\in J}K'^\bullet_j\to L^\bullet$
and $\psi^\bullet:K'^\bullet\to K^\bullet$, and we come down to
checking that $\psi^\bullet$ is a quasi-isomorphism. However,
we have $H^rK^\bullet=\prod_{j\in J}H^rK^\bullet_j$ and likewise
for $H^rK'^\bullet$, for every $r\in\Z$, and under this
identification the map $H^r\psi^\bullet$ equals
$\prod_{j\in J}H^r\psi^\bullet_j$, whence the contention.
A similar argument proves the assertion for direct sums.

(ii): We have commutative diagrams in $\sC(A)$ :
$$
{\diagram
&  K^\bullet_i \ar[r]^-{\psi^\bullet_i} \ar[d] &
K^\bullet_i\oplus K^\bullet_{i+1} \ar[d] \\
0 \ar[r] & \bigoplus_{n\in\N}K^\bullet_n \ar[r]^-{\psi^\bullet} &
\bigoplus_{n\in\N}K^\bullet_n \ar[r] & L^\bullet
\ar[r] & 0
\enddiagram}
\qquad
\text{for every $i\in\N$}
$$
whose bottom rows are the same short exact sequence for
every $i\in\N$, and with $\psi^r_i(x):=(x,-\phi^r_i(x))$
for every $i\in\N$, every $r\in\Z$ and every $x\in K^r_i$.
For every object $C^\bullet$ of $\sD(A)$ set
$$
P(k):=\prod_{n\in\N}\Hom_{\sD(A)}(K^\bullet_n,C^\bullet[k])
\qquad
\tilde\psi(k):=\Hom_{\sD(A)}(\psi^\bullet,C^\bullet[k])
\qquad
\text{for every $k\in\Z$}.
$$
By (i) and remark \ref{rem_long-exact-derHom}, there follows
a long exact sequence :
$$
\Hom_{\sD(A)}(L^\bullet,C^\bullet[k])\to P(k)
\xrightarrow{\ \tilde\psi(k)\ }P(k)\to
\Hom_{\sD(A)}(L^\bullet,C^\bullet[k+1])\to\cdots
$$
and commutative diagrams for every $i\in\N$ and $k\in\Z$ :
\set\begin{equation}\label{eq_as-in-here}
{\diagram P(k) \ar[r]^-{\tilde\psi(k)} \ar[d] & P(k) \ar[d] \\
\Hom_{\sD(A)}(K^\bullet_i\oplus K^\bullet_{i+1},C^\bullet[k])
\ar[r] & \Hom_{\sD(A)}(K^\bullet_i,C^\bullet[k])
\enddiagram}
\end{equation}
whose bottom rows are induced by $\psi_i^\bullet$. Thus,
$\Ker\,\tilde\psi(0)$ consists of the systems of
morphisms $(\beta^\bullet_n:K^\bullet_n\to C^\bullet~|~n\in\N)$
such that $\beta^\bullet_{n+1}\circ\phi^\bullet_n=\beta^\bullet_n$
for every $n\in\N$. In other words, $\Ker\,\tilde\psi(0)$
is naturally identified with the set of cocones with basis
$(K^\bullet_n~|~n\in\N)$ and vertex $C^\bullet$ in the category
$\sD(A)$. Moreover, considering \eqref{eq_as-in-here} with
$k=-1$, a simple induction on $i$ shows that $\tilde\psi(-1)$
is surjective (details left to the reader); summing up, we
have a natural identification
$$
\Hom_{\sD(A)}(L^\bullet,C^\bullet)\isom
\Hom_{\sD(A)}(\colim_{n\in\N}j^I(K^\bullet_n),C^\bullet)
$$
whence the assertion.

(iii): We argue similarly, considering the commutative
diagram in $\sC(A)$ :
$$
{\diagram
0 \ar[r] & L'^\bullet \ar[r] &
\prod_{n\in\N}K^\bullet_n \ar[r]^-{\mu^\bullet} \ar[d] &
\prod_{n\in\N}K^\bullet_n \ar[r] \ar[d] & 0 \\
& & K^\bullet_i\oplus K^\bullet_{i+1} \ar[r]^-{\mu^\bullet_i}
& K^\bullet_i \\
\enddiagram}
\qquad
\text{for every $i\in\N$}
$$
whose top row is independent of $i$ and with
$\mu^r_i(x,y):=x-\phi^r_i(y)$ for every $i\in\N$, every
$r\in\Z$ and every $(x,y)\in K^r_i\oplus K^r_{i+1}$. Notice
that the top row is a short exact sequence in $\sC(A)$,
since the system $(K^r_n~|~n\in\N)$ satisfies the
Mittag-Leffler condition for every $r\in\Z$.
For every $k\in\Z$ and any $C^\bullet\in\Ob(\sD(A))$ set
$$
Q(k):=\prod_{n\in\N}\Hom_{\sD(A)}(C^\bullet,K^\bullet_n[k])
\qquad
\tilde\mu(k):=\Hom_{\sD(A)}(C^\bullet,\mu^\bullet[k]).
$$
By (i) and remark \ref{rem_long-exact-derHom}, there
follows a long exact sequence
$$
\Hom_{\sD(A)}(C^\bullet,L'^\bullet[k])\to
Q(k)\xrightarrow{\ \tilde\mu(k)\ }Q(k)
\to\Hom_{\sD(A)}(C^\bullet,L'^\bullet[k+1])\to
$$
as well as commutative diagrams as in \eqref{eq_as-in-here}.
Then, a simple inspection yields a natural identification
of $\Ker\,\tilde\mu(0)$ with the set of cones in $\sD(A)$
with basis $(K^\bullet_n~|~n\in\N)$ and vertex $C^\bullet$.
Lastly, $\Coker\,\tilde\mu(-1)$ is naturally identified with
$\lim^1_{n\in\N}\Hom_{\sD(A)}(C^\bullet,K^\bullet_n[-1])$, whence the
contention.
\end{proof}

\sset\subsubsection{Minimal resolutions}
Let $A$ be a local ring, $k$ its residue field, and :
$$
\cdots\xrightarrow{d_3}L_2\xrightarrow{d_2}
L_1\xrightarrow{d_1}L_0\xrightarrow{\eps} M
$$
a resolution of an $A$-module $M$ of finite type.
We say that $(L_\bullet,d_\bullet,\eps)$ is a
{\em finite-free resolution\/} if each $L_i$ is a
free $A$-module of finite rank. We say that
$(L_\bullet,d_\bullet,\eps)$ is a
{\em minimal free resolution\/} of $M$ if it is a
finite-free resolution, and moreover the induced maps
$k\otimes_AL_i\to k\otimes_A\Img\,d_i$ are isomorphisms
for all $i\in\N$ (where we let $d_0:=\eps$). One verifies
easily that if $A$ is a coherent ring, then every
finitely presented $A$-module admits a minimal resolution.

\sset\subsubsection{}
Let $\underline{L}:=(L_\bullet,d_\bullet,\eps)$ and
$\underline{L}':=(L'_\bullet,d'_\bullet,\eps')$ be two
free resolutions of $M$. A {\em morphism of resolutions\/}
$\underline{L}\to\underline{L}'$ is a map of complexes
$\phi_\bullet:(L_\bullet,d_\bullet)\to(L'_\bullet,d'_\bullet)$
that extends to a commutative diagram
$$
\xymatrix{
L_\bullet \ar[r]^-\eps \ar[d]_{\phi_\bullet} & M \ddouble \\
L'_\bullet \ar[r]^-{\eps'} & M
}$$

\begin{lemma}\label{lem_min-resol}
Let $\underline{L}:=(L_\bullet,d_\bullet,\eps)$
be a minimal free resolution of an $A$-module $M$ of finite type,
$\underline{L}':=(L'_\bullet,d'_\bullet,\eps')$ any other
finite-free resolution, $\phi_\bullet:\underline{L}'\to\underline{L}$
a morphism of resolutions. Then $\phi_\bullet$ is an epimorphism
(in the category of complexes of $A$-modules), $\Ker\,\phi_\bullet$
is a null homotopic complex of free $A$-modules, and there
is an isomorphism of complexes :
$$
L'_\bullet\isom L_\bullet\oplus\Ker\,\phi_\bullet.
$$
\end{lemma}
\begin{proof} Suppose first that $\underline{L}=\underline{L}'$.
We set $d_0:=\eps$, $L_{-1}:=M$, $\phi_{-1}:=\one_M$ and we show
by induction on $n$ that $\phi_n$ is an isomorphism. Indeed,
this holds for $n=-1$ by definition. Suppose that $n\geq 0$
and that the assertion is known for all $j<n$; by a little
diagram chasing (or the five lemma) we deduce that $\phi_{n-1}$
induces an automorphism $\Img\,d_n\isom\Img\,d_n$, therefore
$\phi_n\otimes_A\one_k:k\otimes_AL_n\to k\otimes_AL_n$
is an automorphism (by minimality of $\underline{L}$), so the
same holds for $\phi_n$ ({\em e.g.} by looking at the determinant
of $\phi_n$).

For the general case, by standard arguments we construct
a morphism of resolutions:
$\psi_\bullet:\underline{L}\to\underline{L}'$. By the foregoing
case, $\phi_\bullet\circ\psi_\bullet$ is an automorphism
of $\underline{L}$, so $\phi_\bullet$ is necessarily an
epimorphism, and $\underline{L'}$ decomposes as claimed.
Finally, it is also clear that $\Ker\,\phi_\bullet$ is an
acyclic bounded above complex of free $A$-modules, hence it
is null homotopic.
\end{proof}

\begin{remark}\label{rem_syzy}
(i) Suppose that $A$ is a coherent local ring, and let
$\underline{L}:=(L_\bullet,d_\bullet,\eps)$ and
$\underline{L}':=(L'_\bullet,d'_\bullet,\eps')$ be two
minimal resolutions of the finitely presented $A$-module $M$.
It follows easily from lemma \ref{lem_min-resol} that
$\underline{L}$ and $\underline{L}'$ are isomorphic as
resolutions of $M$.

(ii) Moreover, any two isomorphisms
$\underline{L}\to\underline{L}'$ are homotopic, hence the
rule: $M\mapsto L_\bullet$ extends to a functor
$$
A\Mod_\coh\to\Hot(A\Mod)
$$
from the category of finitely presented $A$-modules to the
homotopy category of complexes of $A$-modules.

(iii) The sequence of $A$-modules $(\Syz^i_AM:=\Img\,d_i~|~i>0)$
is determined uniquely by $M$ (up to non-unique isomorphism).
The graded module $\Syz^\bullet_AM$ is sometimes
called the {\em syzygy\/} of the module $M$. Moreover, if
$\underline{L}''$ is any other finite free resolution of $M$,
then we can choose a morphism of resolutions
$\underline{L}''\to\underline{L}$, which will be a split
epimorphism by lemma \ref{lem_min-resol}, and the submodule
$d_\bullet(L''_\bullet)\subset L_\bullet$ decomposes as a
direct sum of $\Syz^\bullet_AM$ and a free $A$-module of
finite rank.
\end{remark}

\begin{lemma}\label{lem_syzy} Let $A\to B$ a faithfully flat
homomorphism of coherent local rings, $M$ a finitely presented
$A$-module. Then there exists an isomorphism of graded $B$-modules :
$$
B\otimes_A\Syz^\bullet_AM\to\Syz^\bullet_B(B\otimes_AM).
$$
\end{lemma}
\begin{proof} Left to the reader.
\end{proof}

\begin{proposition}\label{prop_resol-exist}
Let $A\to B$ be a flat and essentially finitely presented
local ring homomorphism of local rings, $M$ an $A$-flat
finitely presented $B$-module. Then the $B$-module $M$
admits a minimal free resolution
$$
\Sigma_\bullet \quad :\quad 
\cdots\xrightarrow{\ d_3\ } L_2\xrightarrow{\ d_2\ } L_1
\xrightarrow{\ d_1\ } L_0\xrightarrow{\ d_{-1}\ } L_{-1}:=M.
$$
Moreover, $\Sigma$ is {\em universally $A$-exact}, {\em i.e.}
for every $A$-module $N$, the complex $\Sigma_\bullet\otimes_AN$
is still exact.
\end{proposition}
\begin{proof} To start out, let us notice :

\begin{claim}\label{cl_in-order}
In order to prove the proposition, it suffices to show that,
for every $A$-flat finitely presented $B$-module $N$, and
every $B$-linear surjective map $d:L\to N$, from a free
$B$-module $L$ of finite rank, $\Ker\,d$ is also an $A$-flat
finitely presented $B$-module.
\end{claim}
\begin{pfclaim} Indeed, in that case, we can build inductively
a minimal resolution $\Sigma_\bullet$ of $M$, such that
$N_i:=\Ker(d_i:L_i\to L_{i-1})$ is an $A$-flat finitely
presented $B$-module for every $i\in\N$. Namely, suppose
that a complex $\Sigma^{(i)}_\bullet$ with these properties has
already been constructed, up to degree $i$, and let $\kappa_B$
be the residue field of $B$; by Nakayama's lemma, we may
find a surjection $d_{i+1}:L_{i+1}\to N_i$, where $L_{i+1}$
is a free $B$-module of rank $\dim_{\kappa_B}(N_i\otimes_B\kappa_B)$.
Under the assumption of the claim, the resulting complex
$\Sigma_\bullet^{(i+1)}:(L_{i+1}\to L_i\to L_{i-1}\to\cdots\to M)$
fulfills the sought conditions, up to degree $i+1$.

It is easily seen that the complex $\Sigma_\bullet$ thus
obtained shall be universally $A$-exact.
\end{pfclaim}

Let us write $A$ as the union of the filtered
family $(A_\lambda~|~\lambda\in\Lambda)$ of its noetherian
subalgebras. Say that $B=C_\fp$, for some finitely presented
$A$-algebra $C$, and a prime ideal $\fp\subset C$, and
$M=N_\fp$ for some finitely presented $C$-module $N$. We
may find $\lambda\in\Lambda$, a finitely generated
$A_\lambda$-algebra $C_\lambda$ and a $C_\lambda$-module
$N_\lambda$ such that $C=C_\lambda\otimes_{A_\lambda}A$, and
$N=N_\lambda\otimes_{A_\lambda}A$; for every $\mu\geq\lambda$,
let $C_\mu:=A_\mu\otimes_{A_\lambda}C_\lambda$ and
$N_\mu:=A_\mu\otimes_{A_\lambda}N_\lambda$; also, denote by
$\fp_\mu$ the preimage of $\fp$ in $C_\mu$, and set
$B_\mu:=(C_\mu)_{\fp_\mu}$, $M_\mu:=(N_\mu)_{\fp_\mu}$.
According to \cite[Ch.IV, Cor.11.2.6.1(ii)]{EGAIV-3}, we may assume
that $M_\mu$ is a flat $A_\mu$-module, for every $\mu\geq\lambda$.
Moreover, suppose that $d:L\to M$ is a $B$-linear surjection
from a free $B$-module $L$ of rank $r$; then we may find
$\mu\in\Lambda$ such that $d$ descends to a $B_\mu$-linear
surjection $d_\mu:L_\mu\to M_\mu$ from a free $B_\mu$-module
$L_\mu$ of rank $r$. It follows easily that $K_\mu:=\Ker\,d_\mu$
is a flat $A_\mu$-module, for every $\mu\geq\lambda$, and the
induced map $K_\mu\otimes_{B_\mu}B\to K:=\Ker\,d$ is a surjection,
whose kernel is a quotient of $\Tor_1^{A_\mu}(A,M_\mu)_\fp$;
the latter vanishes, since $M_\mu$ is $A_\mu$-flat. Hence,
$K$ is $A$-flat; furthermore, $K_\mu$ is clearly a finitely
generated $B_\mu$-module, hence $K$ is a finitely presented
$B$-module. Then the proposition follows from claim
\ref{cl_in-order}.
\end{proof}

\subsection{Simplicial objects}\label{sec_simplicial}
In this section, we introduce the simplicial formalism,
which provides the language for the homotopical algebra
of section \ref{sec_homotopy}.

\begin{definition}\label{def_simplicial-cats}
Let $\cC$ be any category, and $k\in\N$ any integer.
\begin{enumerate}
\item
We denote by $\Delta$ the {\em simplicial category}, whose
objects are the finite ordered sets :
$$
[n]:=\{0<1<\cdots<n\}
\qquad
\text{for every $n\in\N$}
$$
and whose morphisms are the non-decreasing functions.
\item
$\Delta$ is a full subcategory of the {\em augmented
simplicial category\/} $\Delta^{\!\wedge}$, whose set of
objects is $\Ob(\Delta)\cup\{\emptyset\}$, with
$\Hom_{\Delta^{\!\wedge}}(\emptyset,[n])$ consisting of the
unique mapping of sets $\emptyset\to[n]$, for every
$n\in\N$. It is convenient to set $[-1]:=\emptyset$.
\item
The {\em augmented $k$-truncated simplicial category}
$\Delta^{\!\wedge}_k$, is the full subcategory of
$\Delta^{\!\wedge}$ whose objects are the elements of
$\Ob(\Delta^{\!\wedge})$ of cardinality $\leq k+1$.
The {\em $k$-truncated simplicial category\/} is the full
subcategory $\Delta_k$ of $\Delta^{\!\wedge}_k$ whose set
of objects is $\Ob(\Delta^\wedge_k)\setminus\{\emptyset\}$.
\item
A {\em simplicial object\/} (resp. an {\em augmented simplicial
object}, resp. a {\em $k$-truncated simplicial object}, resp.
a {\em $k$-truncated augmented simplicial object}) of $\cC$
is a functor $\Delta^o\to\cC$ (resp. $(\Delta^\wedge)^o\to\cC$,
resp. $\Delta_k^o\to\cC$, resp. $(\Delta^\wedge_k)^o$).
The morphisms of simplicial objects of $\cC$ are just the
natural transformations (and likewise for the truncated or
augmented variants). Clearly, these objects form a category,
and we use the notation
$$
\begin{aligned}
s.\cC
& :=\bFun(\Delta^o,\cC)
& \qquad &
& \hat s.\cC
& :=\bFun((\Delta^\wedge)^o,\cC)
\\
s_k.\cC
& :=\bFun(\Delta^o_k,\cC)
& \qquad &
& \hat s_k.\cC
& :=\bFun((\Delta^\wedge_k)^o,\cC).
\end{aligned}
$$
\item
Dually, a {\em cosimplicial object\/} $F^\bullet$ of $\cC$ is
a functor $F:\Delta\to\cC$, or -- which is the same -- a
simplicial object in $\cC^o$. Likewise one defines the
truncated or augmented cosimplicial variants, and we set
$$
\begin{aligned}
c.\cC
& :=\bFun(\Delta,\cC)
& \qquad &
& \hat c.\cC
& :=\bFun(\Delta^\wedge,\cC)
\\
c_k.\cC
& :=\bFun(\Delta_k,\cC)
& \qquad &
& \hat c_k.\cC
& :=\bFun(\Delta^\wedge_k,\cC).
\end{aligned}
$$
\end{enumerate}
\end{definition}

\sset\subsubsection{}\label{subsec_front-back}
Notice the {\em front-to-back\/} involution :
\set\begin{equation}\label{eq_front-back}
\Delta\to\Delta \quad :\quad
\quad (\alpha:[n]\to[m])\mapsto(\alpha^\vee:[n]\to[m])
\quad\text{for every $n,m\in\N$}
\end{equation}
defined as the endofunctor which induces the identity
on $\Ob(\Delta)$, and such that :
$$
\alpha^\vee(i):=m-\alpha(n-i)\qquad
\text{for every $\alpha\in\Hom_\Delta([n],[m])$ and every $i\in[n]$}.
$$
Another construction of interest is the endofunctor
$$
\gamma:\Delta\to\Delta
$$
given by the rule : $[n]\mapsto[n+1]$ for every $n\in\N$,
and which takes any morphism $\alpha:[n]\to[m]$ of $\Delta$,
to the morphism $\gamma(\alpha):[n+1]\to[m+1]$ which is the
unique extension of $\alpha$ such that $\gamma(\alpha)(n+1):=m+1$.
Notice that $\gamma$ restricts to functors
$\gamma_k:\Delta_{k+1}\to\Delta_k$ for every $k\in\N$.

$\bullet$\ \ 
Given a simplicial object $F$ of $\cC$, one gets a cosimplicial
object $F^o$ of $\cC$, by the (obvious) rule : $(F^o)[n]:=F[n]$
for every $n\in\N$, and $F^o(\alpha):=F(\alpha)^o$ for every
morphism $\alpha$ in $\Delta$.

$\bullet$\ \
Moreover, by composing a simplicial (resp. cosimplicial) object
$F$ (resp. $G$) with the involution \eqref{eq_front-back}, one
obtains a simplicial (resp. cosimplicial) object $F^\vee$ (resp.
$G^\vee$). Likewise, given a morphism $\alpha:F_1\to F_2$, the
Godement product $\alpha^\vee:=\alpha*\eqref{eq_front-back}$ is
a morphism $F_1^\vee\to F_2^\vee$.

$\bullet$\ \
For $F$ and $G$ as above, we may also consider the simplicial
(resp. cosimplicial) object $\gamma F:=F\circ\gamma^o$ (resp.
$\gamma G:=G\circ\gamma$), and this definition extends again
to morphisms, by taking Godement products. The object $\gamma F$
(resp. $\gamma G$) is called the {\em path space\/} of $F$
(resp. of $G$). If $F$ is a $(k+1)$-truncated simplicial object,
then we can consider $\gamma_kF:=F\circ\gamma^o_k$, which is
a $k$-truncated simplicial object (and likewise for truncated
cosimplicial objects).

\sset\subsubsection{}\label{subsec_simplicial-object}
There is an obvious fully faithful functor :
$$
\cC\to s.\cC
\quad :\quad
A\mapsto s.A
\qquad
\text{(\ resp.\quad $\cC\to s_k.\cC
                    \quad :\quad
                    A\mapsto s_k.A$ \ )}
$$
that assigns to each object $A$ of $\cC$ the {\em constant
simplicial object\/} $s.A$ (resp. {\em constant truncated
simplicial object\/} $s_k.A$) such that $s.A[n]:=A$ for every
$n\in\N$ (resp. for every $n\leq k$), and $s.A(\alpha):=\one_A$
for every morphism $\alpha$ of $\Delta$ (resp. of $\Delta_k$).
Of course, we have as well augmented variants $\hat s.A$ and
$\hat s_k.A$, and cosimplicial versions $c.A$, $c_k.A$, $\hat c.A$,
$\hat c_k.A$.

Moreover, we have, for every integer $k\in\N$, the
{\em $k$-truncation functor}
$$
s.\trunc_k:s.\cC\to s_k.\cC
\qquad
\text{(\ resp.\quad $\hat s.\trunc_k:\hat s.\cC\to\hat s_k.\cC$\ )}
$$
that assigns to any simplicial (resp. augmented simplicial)
object $F:\Delta^o\to\cC$ (resp. $F:(\Delta^\wedge)^o\to\cC$)
its composition with the inclusion functor $\Delta^o_k\to\Delta^o$
(resp. $(\Delta^\wedge_k)^o\to(\Delta^\wedge)^o$). Again, we
have as well the corresponding cosimplicial versions $c.\trunc_k$
and $\hat c.\trunc_k$. Also, for every $n\in\N$, we have
the functor
$$
\bullet[n]:s.\cA\to\cA
\qquad
A\mapsto A[n].
$$
Lastly, any functor $\phi:\cB\to\cC$ induces functors
$$
s.\phi:s.\cB\to s.\cC
\qquad
s_k.\phi:s_k.\cB\to s_k.\cC
\quad :\quad
F\mapsto\phi\circ F
$$
and there are of course augmented variants $\hat s.\phi$
and $\hat s_k.\phi$, as well as the corresponding cosimplicial
versions.

\begin{example}\label{ex_simplicies}
(i)\ \
For every integer $k\geq -1$, we denote by $\bDelta_k$
the simplicial set represented by $[k]$. Explicitly,
$\bDelta_k$ is given by the rule
$$
\bDelta_k[n]:=\Hom_\Delta([n],[k])
\quad\text{and}\quad
\bDelta_k[\phi]:=\Hom_\Delta(\phi,[k])
$$
for every $n\in\N$ and every morphism $\phi$ in $\Delta$.
For instance, $\bDelta_{-1}$ (resp. $\bDelta_0$) is the
constant simplicial set associated to the empty set
(resp. to the set with one element); this is also the
initial (resp. final) object for $s.\Set$.

(ii)\ \
Any morphism $\phi:[k]\to[n]$ in $\Delta^\wedge$ induces
a morphism which we shall denote
$$
\bDelta_\phi:=\Hom_{\Delta^\wedge}(-,\phi):\bDelta_k\to\bDelta_n
$$
given by left composition with $\phi$ on $\bDelta_k[i]$,
for every $i\in\N$. 
\end{example}

\sset\subsubsection{}\label{subsec_do-faces}
For given $n\in\N$, and every $i=0,\dots,n$, let
$$
\eps_i:[n-1]\to[n]
\qquad
\text{(resp. $\eta_i:[n+1]\to[n]$)}
$$
be the unique injective map in $\Delta^\wedge$ whose image misses
$i$ (resp. the unique surjective map in $\Delta$ with two elements
mapping to $i$). The morphisms $\eps_i$ (resp. $\eta_i$) are
called {\em face maps} (resp. {\em degeneracy maps}). By
direct inspection, one checks that they fulfill the identities :
$$
\begin{aligned}
\eps_j\circ\eps_i & =\eps_i\circ\eps_{j-1} \qquad
\text{if $i<j$} \\
\eta_j\circ\eta_i & =\eta_i\circ\eta_{j+1} \qquad
\text{if $i\leq j$} \\
\eta_j\circ\eps_i & =
\left\{\begin{array}{ll}
       \eps_i\circ\eta_{j-1}       & \text{if $i<j$} \\
       \one                       & \text{if $i=j$ or $i=j+1$} \\
       \eps_{i-1}\circ\eta_j       & \text{if $i>j+1$}.
       \end{array}\right.
\end{aligned}
$$

\begin{example}\label{ex_face-and-deg}
(i)\ \
For instance, notice the identities :
$$
\eps_i^\vee=\eps_{n-i}
\qquad
\eta_i^\vee=\eta_{n-i}
\qquad\text{for every $n\in\N$ and every $i=0,\dots,n$}.
$$

(ii)\ \
For every $r,s,i\in\N$ with $i\leq r$, we set
$$
\eps_{r,i}^s:=\eps_i\circ\cdots\circ\eps_i:[r]\to[r+s].
$$
This is the injective mapping whose image is
$\{0,\dots,i-1,s+i,\dots,r+s\}$; for instance, the
front-to-back dual
$$
\eps^{s\vee}_{r,0}:=\eps_{r+s}\circ\cdots\circ\eps_{r+1}:[r]\to[r+s]
$$
is just the natural inclusion map. We shall also use the
notation
$$
\eps^s_{-1,0}:[-1]\to[s]
\qquad
\text{for every $s\in\N$}
$$
for the unique morphism $\emptyset\to[s]$ in $\Delta^\wedge$;
of course, we have $\eps^{s\vee}_{-1,0}=\eps^s_{-1,0}$ for every
$s\in\N$.
\end{example}

\sset\subsubsection{}\label{subsec_face-and-degeracies}
It is easily seen that every morphism $\alpha:[n]\to[m]$ in
$\Delta$ admits a unique factorization $\alpha=\eps\circ\eta$,
where the monomorphism $\eps$ is uniquely a composition of faces :
$$
\eps=\eps_{i_1}\circ\cdots\circ\eps_{i_s}
\qquad
\text{with $0\leq i_s\leq\cdots\leq i_1\leq m$}
$$
and the epimorphism $\eta$ is uniquely a composition of
degeneracy maps :
$$
\eta=\eta_{j_1}\circ\cdots\circ\eta_{j_t}
\qquad
\text{with $0\leq j_1<\cdots<j_t\leq m$}
$$
(see \cite[Lemma 8.1.2]{We}). It follows that, to give a
simplicial object $A[\bullet]$ of a category $\cC$, it suffices
to give a sequence of objects $(A[n]~|~n\in\N)$ of $\cC$,
together with {\em face operators}
$$
\partial_i:=A[\eps_i]:A[n]\to A[n-1]
\qquad i=0,\dots,n
$$
for every integer $n>0$ and {\em degeneracy operators}
$$
\sigma_i:=A[\eta_i]:A[n]\to A[n+1]
\qquad i=0,\dots,n
$$
for every $n\in\N$, satisfying the following
{\em simplicial identities\/} :
\set\begin{equation}\label{eq_simpl-identities}
\begin{aligned}
\partial_i\circ\partial_j & =\partial_{j-1}\circ\partial_i \qquad
\text{if $i<j$} \\
\sigma_i\circ\sigma_j & =\sigma_{j+1}\circ\sigma_i \qquad
\text{if $i\leq j$} \\
\partial_i\circ\sigma_j & =
\left\{\begin{array}{ll}
       \sigma_{j-1}\circ\partial_i & \text{if $i<j$} \\
       \one                        & \text{if $i=j$ or $i=j+1$} \\
       \sigma_j\circ\partial_{i-1} & \text{if $i>j+1$}.
       \end{array}\right.
\end{aligned}
\end{equation}
Under this correspondence we have $\partial_i=A(\eps_i)$ and
$\sigma_i=A(\eta_i)$ (\cite[Prop.8.1.3]{We}). Likewise, a
$k$-truncated simplicial object of $\cC$ is the same as the
datum of a sequence $(A[n]~|~n=0,\dots,k)$ of objects of
$\cC$, and of a system of face and degeneracy operators
restricted to this sequence of objects, and fulfilling the
same identities \eqref{eq_simpl-identities}.

Dually, a cosimplicial object $A[\bullet]$ of $\cC$ is the same
as the datum of a sequence $(A[n]~|~n\in\N)$ of objects of $\cC$,
together with {\em coface operators}
$$
\partial^i:A[n-1]\to A[n] \qquad i=0,\dots,n
$$
and {\em codegeneracy operators}
$$
\sigma^i:A[n+1]\to A[n] \qquad i=0,\dots,n
$$
which satisfy the {\em cosimplicial identities} :
\set\begin{equation}\label{eq_cosimpl-identities}
\begin{aligned}
\partial^j\circ\partial^i & =\partial^i\circ\partial^{j-1} \qquad
\text{if $i<j$} \\
\sigma^j\circ\sigma^i & =\sigma^i\circ\sigma^{j+1} \qquad
\text{if $i\leq j$} \\
\sigma^j\circ\partial^i & =
\left\{\begin{array}{ll}
       \partial^i\circ\sigma^{j-1} & \text{if $i<j$} \\
       \one                        & \text{if $i=j$ or $i=j+1$} \\
       \partial^{i-1}\circ\sigma^j & \text{if $i>j+1$}
       \end{array}\right.
\end{aligned}
\end{equation}
and likewise for $k$-truncated cosimplicial objects.

\sset\subsubsection{}\label{subsec_augment-obj}
An augmented simplicial object of a category $\cC$
can be viewed as the datum of a simplicial object $A[\bullet]$
of\/ $\cC$, together with an object $A[-1]\in\Ob(\cC)$, and a
morphism $\eps:A[0]\to A[-1]$, which is an {\em augmentation\/},
{\em i.e.} such that :
$$
\eps\circ\partial_0=\eps\circ\partial_1.
$$
Dually, an augmented cosimplicial object of $\cC$ can be viewed
as a cosimplicial object $A[\bullet]$, together with a morphism
$\eta:A[-1]\to A[0]$ in $\cC$, such that $\eta^o$ is an
augmentation for $A^o[\bullet]$. We say that $\eta$ is an
{\em augmentation} for $A[\bullet]$.

\begin{remark}\label{rem_path-spaces}
Let $A$ be a simplicial object of the category $\cC$.

(i)\ \ 
For every $n\in\N$ we have $\gamma A[n]:=A[n+1]$ (notation of
\eqref{subsec_front-back}), and the face operators
$\gamma A[\eps_i]:\gamma A[n+1]\to\gamma A[n]$ for $i\leq n+1$
(resp. degeneracy operators
$\gamma A[\eta_i]:\gamma A[n]\to\gamma A[n+1]$ for $i\leq n$)
of $\gamma A$ are $\partial_i:A[n+2]\to A[n+1]$ (resp.
$\sigma_i:A[n+1]\to A[n+2]$); {\em i.e.} we drop $\partial_{n+2}$
and $\sigma_{n+1}$. Likewise for the truncated variants.

(ii)\ \
The discarded faces $\partial_{n+2}$ and degeneracies
$\sigma_{n+1}$ can be used to produce natural morphisms
$$
s.A[0]\xrightarrow{\ f_A\ }\gamma A\xrightarrow{\ g_A\ }A.
$$
Namely, we set
$$
f_A[n]:=\sigma_n\circ\cdots\circ\sigma_1
\qquad
g_A[n]:=\partial_{n+1}
\qquad
\text{for every $n\in\N$}.
$$
For every $k\in\N$, the same operation on an object
$A\in\Ob(s_{k+1}.\cC)$ yields natural morphisms
$$
s.\trunc_kA\xrightarrow{\ f_A\ }\gamma_kA\xrightarrow{\ g_A\ }s.\trunc_kA.
$$

(iii)\ \
Since there is a unique morphism $\sigma_{n,0}:[n]\to[0]$
in $\Delta$ for every $n\in\N$, it is easily seen that
the system $(A[\sigma_{n,0}]:A[0]\to A[n])$ defines a
natural morphism
$$
s.A[0]\to A
\qquad
\text{in $s.\cC$}.
$$

(iv)\ \
Likewise, suppose $\eps:A[0]\to A[-1]$ is an augmentation
for $A$; since there is exactly one morphism
$\eps^{n+1}_{-1,0}:\emptyset\to[n]$ in $\Delta^\wedge$ for
every $n\in\N$, we see that the system
$(A[\eps^{n+1}_{-1,0}]~|~n\in\N)$ defines a natural morphism
$$
A\to s.A[-1]
\qquad
\text{in $s.\cC$}.
$$
\end{remark}

\begin{definition}\label{def_simpl-homotopy}
Let $\cC$ be any category. Denote by
$$
e_i:\Delta^o\to[1]/\Delta^o
\qquad
i=0,1
$$
the functor that assigns to each $[n]\in\Ob(\Delta)$ the
unique morphism $[n]\to[1]$ of $\Delta$ whose image is
$\{i\}$. Let also $\st:[1]/\Delta^o\to\Delta^o$ be the
target functor (see \eqref{subsec_slice-cat}).
Let $A$ and $B$ be two simplicial objects of $\cC$, and
$f,g:A\to B$ two morphisms.

\begin{enumerate}
\item
A {\em homotopy} from $f$ to $g$ is the datum of a natural
transformation
$$
u:A\circ\st\Rightarrow B\circ\st
$$
such that $u*e_0=f$ and $u*e_1=g$.
\item
If $A[\bullet]$ is an augmented simplicial object of $\cC$,
with augmentation given by a morphism $\phi:A\to s.A[-1]$
in $s.\cC$, then we say that $A$ is {\em homotopically trivial},
if there exist a morphism $\psi:s.A[-1]\to A$ and homotopies
$u$ and $v$, respectively from $\one_{s.A[-1]}$ to $\phi\circ\psi$,
and from $\one_A$ to $\psi\circ\phi$.
\item
Dually, if $C$ and $D$ are cosimplicial objects of $\cC$, and
$p,q:C\to D$ any two morphisms in $c.\cC$, then a {\em homotopy}
from $p$ to $q$ is a homotopy from $p^o$ to $q^o$ in $s.\cC^o$.
Likewise we define homotopically trivial cosimplicial objects.
\end{enumerate}
\end{definition}

\begin{remark}\label{rem_simpl-homotopies}
(i)\ \
In the situation of definition \ref{def_simpl-homotopy},
suppose that $p:A'\to A$ and $q:B\to B'$ are any two morphisms
in $s.\cC$; then the natural transformation
$$
(q*\st)\circ u\circ(p*\st):A'\circ\st\Rightarrow B'\circ\st
$$
is a homotopy from $f\circ p$ to $q\circ g$. Moreover, if
$F:\cC\to\cD$ is any functor, then
$$
F*u:FA*\st\Rightarrow FB*\st
$$
is a homotopy from $Ff$ to $Fg$.

(ii)\ \
However, unlike the case for chain homotopies, simplicial
homotopies cannot be composed in this generality; hence, the
simplicial category $s.\cC$ cannot be made into a $2$-category,
by taking the homotopies as $2$-cells.

(iii)\ \
In the same vein, for a general category $\cC$, the relation
``there exists a homotopy from $f$ to $g$'' on morphisms of
$s.\cC$ is neither symmetric nor transitive (though it will
follow from theorem \ref{th_Dold-Kan} that this is an equivalence
relation, in case $\cC$ is abelian).
\end{remark}

\sset\subsubsection{}\label{subsec_homotopy-by-us}
Notice that there are exactly $n+1$ morphisms $\phi:[n]\to[1]$
for every $[n]\in\Ob(\Delta)$, and they can be labeled by
the cardinality of $\phi^{-1}(0)$ : for every $n\in\N$ and
every $k\leq n+1$, we shall write $\phi_{n,k}:[n]\to[1]$ for
the unique morphism such that $\phi^{-1}_{n,k}(0)$ has cardinality
$k$. With this notation, notice that
$$
\begin{aligned}
\phi_{n,k}\circ\eps_i=\left\{\begin{array}{ll}
                              \phi_{n-1,k} & \text{if $i\geq k$} \\
                              \phi_{n-1,k-1} & \text{if $i<k$}
                            \end{array}\right.
\end{aligned}
\qquad\text{and}\qquad
\begin{aligned}
\phi_{n,k}\circ\eta_i=\left\{\begin{array}{ll}
                              \phi_{n+1,k}   & \text{if $i\geq k$} \\
                              \phi_{n+1,k+1} & \text{if $i<k$}.
                            \end{array}\right.
\end{aligned}
$$
Hence, a homotopy $u$ from $f$ to $g$ as in definition
\ref{def_simpl-homotopy}, is the same as a system of morphisms
$$
u_{n,k}:A[n]\to B[n]
\qquad
\text{for every $n\in\N$ and every $k\leq n+1$}
$$
such that $u_{n,n+1}=f[n]$ and $u_{n,0}=g[n]$ for every $n\in\N$,
and the diagrams
$$
\xymatrix{ A[n] \ar[rr]^-{u_{n,k}} \ar[d]_{\partial_i} & &
B[n] \ar[d]^{\partial_i} & &
A[n] \ar[rr]^-{u_{n,k}} \ar[d]_{\sigma_i} & & B[n] \ar[d]^{\sigma_i} \\
A[n-1] \ar[rr]^-{u_{n-1,k-a}} & & B[n-1] & &
A[n+1] \ar[rr]^-{u_{n+1,k+a}} & & B[n+1]
}$$
commute for every $n\in\N$ and every $k\leq n+1$, where
$a:=0$ if $i\geq k$, and $a:=1$ if $i<k$.

\sset\subsubsection{}\label{subsec_bisimplex}
A {\em bisimplicial object} in a category $\cC$ is an object
of the category $s.(s.\cC)$. The latter can also be regarded
as the category of all functors $\Delta^o\times\Delta^o\to\cC$;
it follows that a bisimplicial object of $\cC$ is the same as
a system $(A[p,q]~|~(p,q)\in\N\times\N)$ of objects of $\cC$,
together with morphisms
$$
A[\alpha,\beta]:A[p,q]\to A[p',q']
\qquad
\text{for all morphisms
$\alpha:[p']\to[p]$, $\beta:[q']\to[q]$ of $\Delta$}
$$
compatible with compositions of morphisms in $\Delta$, in the
obvious way. More generally, we may define inductively the
category of $n$-simplicial objects $s^n.\cA$, for every $n\in\N$,
by letting $s^0.\cA:=\cA$, and $s^n.\cA:=s.(s^{n-1}.\cA)$, for
every $n>0$. The {\em diagonal functor}
$$
\Delta\to\Delta\times\Delta
\qquad
[n]\mapsto([n],[n])
\qquad
\alpha\mapsto(\alpha,\alpha)
\qquad
\text{for all $n\in\N$ and all morphisms $\alpha$ of $\Delta$}
$$
induces a functor
$$
\sDelta_\cC:s^2.\cC\to s.\cC
\qquad
A\mapsto A^\Delta.
$$
Especially, we have $A^\Delta[n]:=A[n,n]$ for every $n\in\N$, and
the face operators $\partial_i$ on $A^\Delta[n]$ are of the form
$A[\eps_i,\eps_i]$, for every $i=0,\dots,n$ (and likewise
for the degeneracies). Also, the {\em flip functor}
$$
\Delta\times\Delta\to\Delta\times\Delta
\qquad
([m],[n])\mapsto([n],[m])
$$
induces an endofunctor
$$
\mathsf{fl}:s^2.\cA\to s^2.\cA
$$
in the obvious way. Furthermore, the endofunctors
$$
\Delta\times\Delta\xrightarrow{\ \gamma\times\one_\Delta\ }
\Delta\times\Delta\xleftarrow{\ \one_\Delta\times\gamma\ }
\Delta\times\Delta
$$
induce functors
$$
s^2.\cA\xleftarrow{\ \gamma_1\ }s^2.\cA\xrightarrow{\ \gamma_2\ }
s^2.\cA
$$
that admit descriptions as in remark \ref{rem_path-spaces}(i).
Correspondingly, we get natural morphisms
$$
 g^{(i)}_A:\gamma_iA\to A
\qquad
\text{for $i=1,2$ and every $A\in\Ob(s^2.\cA)$}
$$
as in remark \ref{rem_path-spaces}(ii).

\begin{remark}\label{rem_mixed-simpl-tensors}
(i)\ \
Let $\cC$ be a category whose finite coproducts are representable.
Let also $\fSet$ be the category of finite sets. To every
object $S$ of $s.\fSet$ and every $X\in\Ob(s.\cC)$, we attach
a bisimplicial object $S\boxtimes X$ of $\cC$ as follows. For
every $n,m\in\N$, we let $S\boxtimes X[n,m]$ be the coproduct
of finitely many copies of $X[m]$, indexed by the elements of
$S[n]$; hence, for every $a\in S[n]$ we have a natural morphism
$i_a:X[m]\to S\boxtimes X[n,m]$.
If $\phi:[n]\to[n']$ and $\psi:[m]\to[m']$ are any two morphisms
in $\Delta^o$, we let
$S\boxtimes X[\phi,\psi]:S\boxtimes X[n,m]\to S\boxtimes X[n',m']$
be the unique morphism such that
$S\boxtimes X[\phi,\psi]\circ i_a=i_{S[\phi](a)}\circ X[\psi]$
for every $a\in S[n]$. Clearly, this rules extends to a well
defined functor
$$
s.\fSet\times s.\cC\to s^2.\cC
\qquad
(S,X)\mapsto S\boxtimes X.
$$
Likewise, we define $X\boxtimes S:=\mathsf{fl}(S\boxtimes X)$
(notation of \eqref{subsec_bisimplex}).
If all coproducts of $\cC$ are representable, we may even
extend the above construction to arbitrary simplicial sets.

(ii)\ \
In the same vein, let $\cA$ be any abelian category, and
$M$ any object of $s.\Z\Mod_\mathrm{fg}$ (notation of
\eqref{subsec_mixed-tensors}). For any $A\in\Ob(s.\cA)$,
we may define a bisimplicial object $M\boxtimes_\Z A$
of $\cA$, by the rule : $[n,m]\mapsto M[n]\otimes_\Z A[n]$
for every $n,m\in\N$ and
$[\phi,\psi]\mapsto M[\phi]\otimes_\Z A[\psi]$ for all
morphisms $\phi,\psi$ of $\Delta$ (where these mixed tensor
products are as defined in \eqref{subsec_mixed-tensors}).
Clearly these rules yield a well defined functor
$$
s.\Z\Mod_\mathrm{fg}\times s.\cA\to s^2.\cA
\qquad
(M,A)\mapsto M\boxtimes_\Z A.
$$
Likewise, we set $A\boxtimes_\Z M:=\mathsf{fl}(M\boxtimes_\Z A)$.

(iii)\ \
Furthermore, if $(\cC,\otimes)$ is any tensor category,
and $X,Y$ any two simplicial objects of $\cC$, we may define
a bisimplicial object $X\boxtimes Y$, by the same rule as
in (ii). This yields a functor
$$
s.\cC\times s.\cC\to s^2.\cC
\qquad
(X,Y)\mapsto X\boxtimes Y.
$$
In this situation (resp. in the situation of (i), resp. of (ii)),
we shall let also
$$
X\otimes Y:=(X\boxtimes Y)^\Delta
\qquad
\text{(resp.\ \ $S\otimes X:=(S\boxtimes X)^\Delta$,\ \ 
resp.\ \ $M\otimes_\Z A:=(M\boxtimes_\Z A)^\Delta$\ )}
$$
which we shall call the {\em tensor product} of $X$ and $Y$
(resp. of $S$ and $X$, resp. of $M$ and $A$). Notice the
natural identifications
$$
S\boxtimes A\isom(S\otimes s.\Z)\boxtimes_\Z A
\qquad
\text{for every $S\in\Ob(s.\fSet)$ and every $A\in\Ob(s.\cA)$}.
$$
Likewise, if $U$ is any unit object for $(\cC,\otimes)$,
we get natural identifications
$$
S\boxtimes X\isom(S\otimes s.U)\boxtimes X
\qquad
\text{for every $S\in\Ob(s.\fSet)$ and every $X\in\Ob(s.\cC)$}
$$
via the isomorphisms $u_X:X\isom U\otimes X$ provided by
proposition \ref{prop_unit-isom}. In the same vein, we have
a natural identification $\bDelta_0\otimes X\isom X$ for
every $X\in\Ob(\cC)$ (notation of example \ref{ex_simplicies}(i)).

(iv)\ \
Notice also that if $f,g:S\to T$ are any two morphisms
of simplicial finite sets, then -- in light of remark
\ref{rem_simpl-homotopies}(i) -- any homotopy $u$ from $f$
to $g$ induces a homotopy $u\otimes X$ from $f\otimes X$ to
$g\otimes X$. 

(v)\ \
Let $f,g:A\to B$ be as in definition \ref{def_simpl-homotopy}.
Notice that the datum of a homotopy $u$ from $f$ to $g$
is the same as that of a morphism
$$
\tilde u:\bDelta_1\otimes A\to B
\qquad
\text{such that $\tilde u\circ(\bDelta_{\eps_0}\otimes\one_A)=f$
            and $\tilde u\circ(\bDelta_{\eps_1}\otimes\one_A)=g$}
$$
(where $\eps_0,\eps_1:[0]\to[1]$ are the face maps : notation
of example \ref{ex_simplicies}(ii)).
Indeed, given $u$, we construct $\tilde u$ as follows.
For every $n\in\N$ and every $\phi\in\bDelta_1[n]$, let
$\tilde u[n]$ be the unique morphism such that
$\tilde u[n]\circ i_\phi=u_\phi$. The naturality of $u$
easily implies that this rule amounts to a morphism $\tilde u$
as sought. Conversely, given $\tilde u$, we can construct
a natural transformation $u$, by reversing the foregoing rule.
\end{remark}

\begin{remark}\label{rem_cosimpl-simpl}
(i)\ \ 
Let $(\cC,\otimes)$ be a tensor category with internal $\Hom$
functor $\cHom$, and unit object $U$. To every two objects $X,Y$
of $s.\cC$, we may attach the object of $\cC$
$$
\cHom_{s.\cC}(X,Y):=\xymatrix{\Equal(\prod_{n\in\N}\cHom(X[n],Y[n])
\ar@<.5ex>[r]^-{d_1} \ar@<-.5ex>[r]_-{d_0} &
\prod_{\phi:[n]\to[m]}\cHom(X[m],Y[n]) )}
$$
where the second product ranges over the morphisms $\phi$ of
$\Delta$, and where
$$
d_0:=\prod_{\phi:[n]\to[m]}\cHom(X[\phi],Y[n])
\qquad\text{and}\qquad
d_1:=\prod_{\phi:[n]\to[m]}\cHom(X[m],Y[\phi]).
$$
Arguing as in example \ref{ex_monoidal}(vi), it is easily seen
that there are natural isomorphisms
$$
\Hom_\cC(Z,\cHom_{s.\cC}(X,Y))\isom\Hom_{s.\cC}(s.Z\otimes X,Y)
$$
for every $Z\in\Ob(\cC)$ and every $X,Y\in\Ob(s.\cC)$.

(ii)\ \
Suppose additionally, that all finite coproducts of $\cC$ are
representable. For every $Z\in\Ob(\cC)$, and every $X\in\Ob(s.\cC)$,
consider the simplicial set
$$
\Hom_\cC(Z,X)
$$
such that $\Hom_\cC(Z,X)[n]:=\Hom_\cC(Z,X[n])$ for every $n\in\N$,
with face and degeneracies deduced from those of $X$, in the obvious
way. For any $k\in\N$, let also $\bDelta_k$ be the simplicial set
defined in example \ref{ex_simplicies}(i); we have natural isomorphisms
$$
\begin{aligned}
\Hom_\cC(Z,\cHom_{s.\cC}(\bDelta_k\otimes s.U,X))
\isom\, & \Hom_{s.\cC}(\bDelta_k\otimes s.Z,X) \\
\isom\, &
\Hom_{s.\Set}(\bDelta_k,\Hom_\cC(Z,X)) \\
\isom\, & \Hom_\cC(Z,X)[k]
\end{aligned}
$$
where the last isomorphism follows from Yoneda's lemma
(proposition \ref{prop_yoneda}(ii) : details left to the reader).
Applying again Yoneda's lemma, we deduce a natural isomorphism
$$
\cHom_{s.\cC}(\bDelta_k\otimes s.U,X)\isom X[k]
\qquad
\text{for every $k\in\N$}.
$$

(iii)\ \
Notice that $\bDelta_\psi\circ\bDelta_\phi=\bDelta_{\psi\circ\phi}$
for every two morphisms $\phi:[k]\to[k']$ and $\psi:[k']\to[k'']$
of $\Delta^\wedge$ (notation of example \ref{ex_simplicies}(ii)).
Hence, the system $(\bDelta_i~|~i\in\N)$ amounts to an
object of $c.s.\Set$. Moreover, since the Yoneda isomorphisms
are natural in both $X$ and $\bDelta_k$, we get a commutative
diagram
$$
\xymatrix{
\cHom_{s.\cC}(\bDelta_{k'}\otimes s.U,X) \ar[r]^-\sim
\ar[d]_{\cHom_{s.\cC}(\bDelta_\phi\otimes s.U,X)} &
X[k'] \ar[d]^{X[\phi]} \\
\cHom_{s.\cC}(\bDelta_k\otimes s.U,X) \ar[r]^-\sim & X[k]
}$$
for every morphism $\phi$ as above.

(iv)\ \
Notice as well that the considerations of (ii) and (iii) can
be repeated, {\em mutatis mutandis}, for truncated simplicial
objects : if $X$ and $Y$ are objects of $s_n.\cC$ and $Z$
is an object of $\cC$, then $\cHom_{s_n.\cC}(X,Y)$
(resp. $\Hom_\cC(Z,X)$) shall be an object of $\cC$ (resp. of
$s_k.\Set$), and we shall have natural isomorphisms
$$
\cHom_{s_n.\cC}(s.\trunc_n(\bDelta_k\otimes s.U),X)\isom X[k]
\qquad
\text{for every $k\leq n$}
$$
and similarly for the commutative diagrams of (iii) (details
left to the reader).
\end{remark}

\sset\subsubsection{}\label{subsec_coskeleton}
Let $\cC$ be a finitely cocomplete category. Theorem
\ref{th_Kan-ext} say that, for every integer $k\in\N$,
the $k$-truncation functor on $s.\cC$ admits a left adjoint
$$
\sk_k:s_k.\cC\to s.\cC
$$
which is called the {\em $k$-th skeleton functor}. By inspecting
the proof of {\em loc.cit.} we see that, for every $k$-truncated
simplicial object $F$, this adjoint is calculated by the rule :
\set\begin{equation}\label{eq_sk}
\sk_kA[n]:=\colim_{\phi:[i]\to[n]}F[i]
\end{equation}
where $i$ ranges over all the integers $\leq k$, and
$\phi$ over all the morphisms $[i]\to[n]$ in $\Delta^o$,
and the transition maps $F[i]\to F[j]$ in the colimit are
the morphisms $F[\psi]$ given by all commutative triangles
\set\begin{equation}\label{eq_comu-triangle}
{\diagram
[i] \ar[rr]^-\psi \ar[rd]_\phi & & [j] \ar[ld]^\phi \\
& [n]
\enddiagram}
\qquad\text{in $\Delta^o$}.
\end{equation}
A morphism $\alpha:[n]\to[m]$ in $\Delta^o$ induces a
morphism $\sk_kF[\alpha]:\sk_kA[n]\to\sk_kA[m]$; namely,
for every $\phi:[i]\to[m]$ one has a natural morphism
$j_\phi:F[i]\to\sk_kF[m]$, and $\sk_kF[\alpha]$
is the colimit of the system of morphisms
$$
j_{\alpha\circ\psi}:F[i]\to\sk_kA[m]
\qquad
\text{for all $\psi:[i]\to[n]$}.
$$
It is clear that, for every $n,m\leq k$ and every
$\alpha:[n]\to[m]$, the colimit \eqref{eq_sk} is realized
by $F[n]$, and under this identification, $\sk_kA[\alpha]$
agrees with $F[\alpha]$, so the unit of adjunction
$$
F\to s.\trunc_k\circ\sk_kF
$$
is an isomorphism. Dually, if all finite limits are
representable in $\cC$, the truncation functor admits
a right adjoint
$$
\cosk_k:s_k.\cC\to s.\cC
$$
called the {\em $k$-th coskeleton functor}, and a simple
inspection of the proof of {\em loc.cit.} yields the rule:
$$
\cosk_kF[n]:=\lim_{\phi:[n]\to[i]}F[i]
$$
where $i$ ranges over the integers $\leq k$, and $\phi:[n]\to[i]$
over the morphisms in $\Delta^o$, and the transition maps are
as in the foregoing (except that the downwards arrows in the
commutative triangles \eqref{eq_comu-triangle} are reversed).
Especially, we easily deduce that the counit of adjunction
$$
s.\trunc_k\circ\cosk_kF\to F
$$
is an isomorphism. Moreover, if $\cB$ is another category with
small $\Hom$-sets, whose finite colimits (resp. finite limits)
are all representable, and $\phi:\cB\to\cC$ is any functor,
then for any $F\in\Ob(s_k.\cB)$ there is a natural transformation
$$
\sk_k(s_k.\phi F)\to s.\phi(\sk_k F)
\qquad
\text{(\ resp.\ $s.\phi(\cosk_k F)\to\cosk_k(s_k.\phi F)$\ )}
$$
which is an isomorphism, if $\phi$ is right exact (resp.
if $\phi$ is left exact).

\sset\subsubsection{}\label{subsec_check-unnorm-diff}
Let $\cA$ be an abelian category, and $A$ any object of
$s.\cA$. For every $n>0$, set
$$
d_n:=\sum^n_{i=0}(-1)^i\cdot\partial_i
\qquad
A[n]\to A[n-1].
$$
Directly from the simplicial identities \eqref{eq_simpl-identities}
we may compute
$$
\begin{aligned}
d_n\circ d_{n+1}=\, &
\sum_{i=0}^n\sum_{j=0}^{n+1}(-1)^{i+j}\cdot\partial_i\circ\partial_j \\
=\, &
\sum_{i=0}^n\sum_{j>i}^{n+1}(-1)^{i+j}\cdot\partial_{j-1}\circ\partial_i
+\sum_{i=0}^n\partial_i\circ\partial_i
+\sum_{i=0}^n\sum_{j<i}^n(-1)^{i+j}\cdot\partial_i\circ\partial_j \\
=\, &
\sum_{i=0}^n\sum_{j-1>i}^{n+1}(-1)^{i+j}\cdot\partial_{j-1}\circ\partial_i
+\sum_{i=0}^n\sum_{j<i}^n(-1)^{i+j}\cdot\partial_i\circ\partial_j \\
=\, & 0
\end{aligned}
$$
for every $n\in\N$, so we are led to the following :

\begin{definition}\label{def_Dold-Kan}
Let $\cA$ be an abelian category, and $k\in\N$ any integer.
\begin{enumerate}
\item
If $A[\bullet]$ is any simplicial object of $\cA$, with
face operators $\partial_i$, the
{\em unnormalized complex associated to $A[\bullet]$}
is the complex $(A_\bullet,d_\bullet)\in\Ob(\sC^{\leq 0}(\cA))$
such that $A_n:=A[n]$ for every $n\in\N$, and $d_n$ is defined
as in \eqref{subsec_check-unnorm-diff}, for every $n>0$.
\item
If $A[\bullet]$ is a $k$-truncated simplicial object of
$\cA$, the {\em unnormalized complex associated to
$A[\bullet]$} as the complex $(A_\bullet,d_\bullet)$ where
$A_n$ and $d_n$ are defined as in (i) for every $n\leq k$,
and $A_n:=0$, $d_n:=0$ for every $n>k$.
\item
If $A[\bullet]$ is as in (i) (resp. as in (ii)), the
{\em normalized complex associated to $A[\bullet]$}
is the subcomplex $N_\bullet A$ of $A_\bullet$ such that
$$
N_0A:=A[0]
\quad\text{and}\quad
N_nA:=\bigcap_{i=1}^n\Ker\,\partial_i
\qquad
\text{for every $n>0$}
\qquad
\text{(resp. for $0<n\leq k$)}.
$$
So, the differential $N_nA\to N_{n-1}A$ equals $\partial_0$,
for every $n\in\N$ (resp. for every $n\leq k$).
\item
If $A[\bullet]$ is as in (i) or (ii), the {\em homology of $A$}
in degree $n$ is
$$
H_nA:=H_nA_\bullet
\qquad
\text{for every $n\in\N$}.
$$
\item
Let $\eps:A\to A_{-1}$ be an augmentation for $A[\bullet]$.
One says that the augmented simplicial object $(A,\eps)$ is
{\em aspherical}, if $H_nA=0$ for every $n>0$, and $\eps$
induces an isomorphism $H_0A\isom A_{-1}$.
\end{enumerate}
\end{definition}

\sset\subsubsection{}\label{subsec_shuffle}
Let $\cA$ be an abelian category, and recall that $s^n.\cA$
and $s_k.\cA$ are both abelian categories as well, for every
$n,k\in\N$ (remark \ref{rem_Add.Fun}(ii)); also, clearly the
rule of definition \ref{def_Dold-Kan}(iii) yields natural
additive functors
$$
\sN_\cA:s.\cA\to\sC^{\leq 0}(\cA)
\qquad
\sN_{\cA,k}:s_k.\cA\to\sC^{[-k,0]}(\cA)
\qquad
A[\bullet]\mapsto N_\bullet A
\qquad
\text{for every $k\in\N$}
$$
and the rules of definition \ref{def_Dold-Kan}(i,ii) yield
additive functors
$$
\sU_\cA:s.\cA\to\sC^{\leq 0}(\cA)
\qquad
\sU_{\cA,k}:s_k.\cA\to\sC^{[-k,0]}(\cA)
\qquad
A[\bullet]\mapsto A_\bullet
\qquad
\text{for every $k\in\N$}.
$$

\begin{remark}\label{rem_degenerate-sub}
Let $\cA$ and $A$ be as in \eqref{subsec_check-unnorm-diff};
directly from the simplicial identities \eqref{eq_simpl-identities}
we may compute
$$
d_n\circ\sigma_k=\sum^{k-1}_{i=0}
(-1)^i\cdot\sigma_{k-1}\circ\partial_i+
\sum^n_{i=k+2}(-1)^i\cdot\sigma_k\circ\partial_{i-1}.
$$
Especially, if we let
$$
D_0:=0
\qquad
\text{and}
\qquad
D_nA:=\sum^{n-1}_{i=0}\Img\,(\sigma_i:A[n-1]\to A[n])
\qquad
\text{for every $n>0$}
$$
we see that $d_n$ restricts to a morphism $d'_n:D_nA\to D_{n-1}A$
for every $n>0$, hence $(D_\bullet A,d'_\bullet)$ is a subcomplex
of $A_\bullet$, called the {\em degenerate subcomplex}, and
clearly we obtain an additive functor
$$
\sD_\cA:s.\cA\to\sC(\cA)
\qquad
A\mapsto D_\bullet A.
$$
\end{remark}

\begin{proposition}\label{prop_degen-plus-norm}
The natural injections induce a decomposition
$$
A_\bullet=N_\bullet A\oplus D_\bullet A
\qquad
\text{in $\sC^{\leq 0}(\cA)$}.
$$
\end{proposition}
\begin{proof} First, we notice that $N_nA\cap D_nA=0$
for every $n\in\N$. Indeed, it suffices to check that
$N_nA\cap\Img(\sigma_i)=0$ for every $i=0,\dots,n-1$,
but the latter follows easily from the identity
$\partial_{i+1}\circ\sigma_i=\one_{A[n-1]}$. To conclude
the proof, it suffices to show that $N_nA+D_nA=A_n$ for
every $n\in\N$. Indeed, set $K_0:=A_n$ and define inductively
$K_i:=K_{i-1}\cap\Ker\,\partial_i$ for every $i=1,\dots,n$;
to prove the latter assertion, it suffices to check that
\set\begin{equation}\label{eq_kernels}
(\one_{A_n}-\sigma_i\cdot\partial_i)(K_i)\subset K_{i+1}
\qquad
\text{for every $i=0,\dots,n$}.
\end{equation}
However, we have :
$$
\partial_j\circ(\one_{A_n}-\sigma_i\cdot\partial_i)=
\partial_j-\sigma_{i-1}\circ\partial_{i-1}\circ\partial_j
\qquad
\text{whenever $j<i$}
$$
and $\partial_i\circ(\one_{A_n}-\sigma_i\cdot\partial_i)=
\partial_i-\partial_i=0$, whence \eqref{eq_kernels}.
\end{proof}

\begin{example}\label{ex_simplicial-exponent}
(i)\ \
Let $S$ be any simplicial set, and set $\Z^S:=S\otimes s.\Z$
(notation of remark \ref{rem_mixed-simpl-tensors}(iii)); we
wish to give an explicit description of the complex
$N_\bullet\Z^S$ . To begin with, notice that $\Z^S[n]$ is the
free abelian group with basis indexed by $S[n]$, for every
$n\in\N$, and for every $x\in S[n]$ denote by $e_x\in\Z^S[n]$
the corresponding basis element. Next, set
$$
D_nS:=\bigcup_{i=0}^{n-1}\Img\,(\sigma_i:S[n-1]\to S[n])
\quad\text{and}\quad
N_nS:=S[n]\setminus D_nS
\qquad
\text{for every $n\in\N$}.
$$
From proposition \ref{prop_degen-plus-norm}, we see that
$N_n\Z^S$ can be naturally identified with the direct
summand of $\Z^S[n]$ generated by the system
$(e_x~|~x\in N_nS)$, for every $n\in\N$. In order to
describe the differential $d_n$, let us define
$\bar e_x:=e_x$ if $x\in N_nS$, and $\bar e_x:=0$ if
$x\in D_nS$. Then we may write
$$
d_n(\bar e_x)=\sum^n_{i=0}(-1)^i\cdot\bar e_{\partial_ix}
\qquad
\text{for every $n\in\N$ and every $x\in N_nS$}.
$$
The verifications are straightforward, and shall be left
to the reader.

(ii)\ \
For instance, for any $i\in\N$ consider the simplicial set
$\bDelta_i$ of example \ref{ex_simplicies}(i), and set
$$
\K\La i\Ra_\bullet:=N_\bullet\Z^{\bDelta_i}
$$
where $\Z^{\bDelta_i}$ is defined as in (i). Notice that this
notation agrees with that of remark \ref{rem_chain-homotopies}(ii).
It is easily seen that a morphism $[n]\to[i]$ of $\Delta$ lies in
$D_n\bDelta_i$ if and only if it is not injective; hence
$N_n\bDelta_i$ is the set of all injective maps $\phi:[n]\to[i]$
of $\Delta$. We deduce a natural isomorphism
$$
\K\La i\Ra_n\isom\Lambda^{n+1}_\Z\Z^{\oplus i+1}
\qquad
\text{for every $i,n\in\N$}.
$$
Namely, to any map $\phi$ as above, we assign the exterior
product $e_{\phi(0)}\wedge\cdots\wedge e_{\phi(n)}$, where
$e_0,\dots,e_i$ denotes the canonical basis of $\Z^{\oplus i+1}$.
Under this isomorphism, the differential of $\K\La i\Ra_\bullet$
gets identified with the differential of the Koszul complex
attached to the sequence $\one_{i+1}:=(1,\dots,1)\in\Z^{\oplus i+1}$,
that will be introduced in \eqref{sec_koszul-alg}. Summing up,
we obtain natural short exact sequences
$$
0\to\K\La i\Ra_\bullet[1]\to\bK_\bullet(\one_{i+1})\to\Z[0]\to 0
\qquad
\text{for every $i\in\N$}
$$
where $[1]$ denotes the shift operator, and $\Z[0]$ is the
complex with $\Z$ in degree zero : see \eqref{subsec_shift}.

(iii)\ \
Let $\cA$ be any abelian category, $A\in\Ob(s.\cA)$ and
$Z\in\Ob(\cA)$. Notice that the simplicial set $\Hom_\cA(Z,A)$
of remark \ref{rem_cosimpl-simpl}(ii) is actually a
simplicial abelian group, and a direct inspection of the
arguments of {\em loc.cit.} yields a natural isomorphism
of abelian groups
$$
\Hom_{s.\Z\Mod}(\Z^{\bDelta_i},\Hom_\cA(Z,A))\isom\Hom_\cA(Z,A[i])
\qquad
\text{for every $i\in\N$}.
$$
\end{example}

\sset\subsubsection{}\label{subsec_sj-and-sq}
Keep the notation of remark \ref{rem_degenerate-sub}, and let
$$
\sj^A_\bullet:N_\bullet A\to A_\bullet
\qquad
\sq^A_\bullet:A_\bullet\to N_\bullet A
$$
be respectively the injection and the projection with kernel
$D_\bullet A$. Then the rules : $A\mapsto\sj^A_\bullet$ and
$A\mapsto\sq^A_\bullet$ define natural transformations
$$
\sj_\bullet:\sN_\cA\Rightarrow\sU_\cA
\qquad
\sq_\bullet:\sU_\cA\Rightarrow\sN_\cA.
$$
With this notation, we may state :

\begin{theorem}\label{th_degen-plus-norm}
With the notation of remark \eqref{subsec_sj-and-sq}, we have :
\begin{enumerate}
\item
The injection $\sj^A_\bullet$ is a homotopy equivalence, and
$D_\bullet A$ is homotopically trivial.
\item
More precisely, there exist natural modifications
$$
\sj_\bullet\circ\sq_\bullet\leadsto\one_{\sU_\cA}
\qquad
\sq_\bullet\circ\sj_\bullet\leadsto\one_{\sN_\cA}
$$
on the $2$-category $\sC(\cA)$ as in example
{\em\ref{ex_Hom-complex-hots}(ii)} (see definition
{\em\ref{def_modification}}).
\item
Especially, the natural map
$$
H_iN_\bullet A\to H_iA
$$
is an isomorphism, for every $i\in\N$.
\end{enumerate}
\end{theorem}
\begin{proof} Clearly (ii)$\Rightarrow$(i), and (i)$\Rightarrow$(iii),
by virtue of remark \ref{rem_abel-homotopies}(ii). To show (ii), set
$\phi^A_0:=\one_{A_0}$ and
$\phi^A_n:=\one_{A[n]}-\sigma_{n-1}\circ\partial_n:A_n\to A_n$
for every $n>0$. We notice :

\begin{claim}\label{cl_supply-hot}
The system $(\phi^A_n~|~n\in\N)$ defines an endomorphism of
$A_\bullet$, which is homotopically equivalent to $\one_{A_\bullet}$.
\end{claim}
\begin{pfclaim} Indeed, let $s_n:=(-1)^n\cdot\sigma_n:A_n\to A_{n+1}$
for every $n\in\N$, and define $s_n:=0$ for every $n<0$. Using
the simplicial identities \eqref{eq_simpl-identities} we compute :
$$
\begin{aligned}
s_{n-1}\circ d_n+d_{n+1}\circ s_n=\, &
(-1)^{n-1}\cdot\sum^n_{j=0}(-1)^j\cdot
(\sigma_{n-1}\circ\partial_j-\partial_j\circ\sigma_n) \\
=\, &
\phi^A_n-\one_{A[n]}+(-1)^{n-1}\cdot\sum^{n-1}_{j=0}(-1)^j\cdot
(\sigma_{n-1}\circ\partial_j-\sigma_{n-1}\circ\partial_j) \\
=\, & \phi^A_n-\one_{A[n]}
\end{aligned}
$$
for every $n>0$, and for $n=0$ we have as well
$s_{-1}\circ d_0+d_1\circ s_0=0$, whence the claim.
\end{pfclaim}

Notice that the homotopy exhibited in the proof of claim
\ref{cl_supply-hot} is natural in $A$, so it already
yields the first sought modification.
Next, for every simplicial object $B$ of $\cA$, let
$\gamma B$ and $g_B:\gamma B\to B$ be as in remark
\ref{rem_path-spaces}(i,ii), and define inductively
$$
\beta^0A:=A
\qquad
\text{and}
\qquad
\beta^{n+1}A:=\Ker\,g_{\beta^nA}
\qquad
\text{for every $n\in\N$}.
$$
A simple inspection shows that
$$
N_{n+i}A\subset \beta^nA[i]
\qquad\text{and}\qquad
\beta^nA[0]=N_nA
\qquad
\text{for every $i,n\in\N$}.
$$
Hence, we may define a subcomplex $B^{(n)}_\bullet$ of $A_\bullet$
for every $n\in\N$, as follows. For $i=0,\dots,n$ we let
$B^{(n)}_i:=N_iA$, and for $i>n$ we let $B^{(n)}_i:=(\beta^nA)[n+i]$
(notation of \eqref{subsec_shift}).
The differential of $B^{(n)}_\bullet$ is of course just the
restriction of that of $A_\bullet$. By construction, $N_\bullet A$
is a subcomplex of $B^{(n)}_\bullet$, for all $n\in\N$. Next, we
define an endomorphism $f^{(n)}_\bullet$ of $B^{(n)}_\bullet$, by
setting
$$
f^{(n)}_i:=\one_{B^{(n)}_i}
\qquad
\text{for every $i\leq n$, and}
\qquad
f^{(n)}_i:=\phi^{\beta^nA}_i
\qquad
\text{for $i>n$}.
$$
Notice that the restriction of $f^{(n)}_\bullet$ to the subcomplex
$N_\bullet A$ is just the inclusion map $N_\bullet A\to B^{(n)}_\bullet$.
Clearly $B^{(n+1)}\subset B^{(n)}$ for every $n\in\N$, and we
remark that $f^{(n)}_\bullet$ factors through the inclusion
map $B^{(n+1)}\subset B^{(n)}$. Indeed, since $B^{(n)}_i=B^{(n+1)}_i$
for every $i\leq n$, the assertion is obvious for this range
of degrees; so we have only to check that $\phi^{\beta^nA}_\bullet$
factors through $\beta^{n+1}A_\bullet$ for every $n\in\N$, and an
easy induction reduces to checking that $\phi^A_\bullet$ factors
through $B^{(1)}_\bullet$. But the latter assertion comes down
to the identity $\partial_i\circ\phi_i^A=0$ for every $i>0$,
which follows easily from the simplicial identities
\eqref{eq_simpl-identities}.

If $(s_i~|~i\in\N)$ is the homotopy between $\one_{\beta^nA_\bullet}$
and $\phi^{\beta^nA}_\bullet$ supplied by claim \ref{cl_supply-hot},
then we obtain a homotopy $(t^{(n)}_i~|~i\in\N)$ between
$\one_{B^{(n)}_\bullet}$ and $f^{(n)}_\bullet$, by setting
$$
t^{(n)}_i:=0
\qquad
\text{for every $i<n$, and}\qquad
t^{(n)}_i:=s_{i-n}
\qquad
\text{for $i\geq n$}.
$$
Notice that $\Img\,t^{(n)}_i\subset D_{i+1}A$ for every $i,n\in\N$.
Thus, for every $n\in\N$, let
$h^{(n)}_\bullet:B^{(n)}_\bullet\to B^{(n+1)}_\bullet$ be the
morphism of complexes deduced from $f^{(n)}_\bullet$ and
$j^{(n)}_\bullet:B^{(n)}\to A_\bullet$ the inclusion map, and set
$q^{(n)}_\bullet:=h^{(n)}_\bullet\circ\cdots\circ h^{(0)}_\bullet$;
it follows easily that the composition
$$
p^{(n)}_\bullet:=
j^{(n+1)}_\bullet\circ q^{(n)}_\bullet:A_\bullet\to A_\bullet
$$
is homotopically equivalent to $\one_{A_\bullet}$, for every $n\in\N$.
More precisely, a direct inspection shows that the system
of morphisms
$$
\tau^{(n)}_i:=
\sum_{k=0}^{n-1}j^{(k+1)}_\bullet\circ t^{(k+1)}_\bullet\circ q^{(k)}_\bullet
\qquad
\text{for every $i\in\N$}
$$
provides a homotopy between $\one_{A_\bullet}$ and $p^{(n)}_\bullet$.
Furthermore, it is clear that
$$
p^{(n)}_i=p^{(m)}_i
\qquad\text{and}\qquad
\tau^{(n)}_i=\tau^{(m)}_i
\qquad
\text{for every $m\geq n$ and every $i\leq n$}
$$
so, we finally get an endomorphism $p_\bullet:A_\bullet\to A_\bullet$
by setting $p_i:=p_i^{(i)}$ for every $i\in\N$, and a homotopy
$\tau_\bullet$ between $p_\bullet$ and $\one_{A_\bullet}$, with
$\tau_i:=\tau_i^{(i)}$ for every $i\in\N$. By construction,
$p_\bullet$ factors through $\sj^A_\bullet$, and moreover, the restriction
of $p_\bullet$ to the subcomplex $N_\bullet A$ is just $\sj^A_\bullet$,
so $\sj^A_\bullet$ is a homotopy equivalence, via natural homotopies,
as stated.

Lastly, notice that $\Img\,\tau_i\subset D_{i+1}A$ for every
$i\in\N$, from which it follows easily that $p_i(D_iA)\subset D_iA$
for every $i\in\N$, and then the foregoing implies that
$\Ker\,p_\bullet=D_\bullet A$. We conclude that $D_\bullet A$
is homotopically trivial, as stated.
\end{proof}

\sset\subsubsection{}\label{subsec_iterate-U}
By iterating $\sU$, we get a functor from
bisimplicial objects to double complexes
$$
\sU^2_\cA:s^2.\cA\xrightarrow{\ s.\sU_\cA\ }s.\sC(\cA)
\xrightarrow{\ \sU_{\sC(\cA)}\ }\sC(\sC(\cA))
\qquad
A[\bullet,\bullet]\mapsto A_{\bullet\bullet}
$$
and notice that this functor is naturally isomorphic to
the functor
$$
s^2.\cA\xrightarrow{\ \sU_{s.\cA}\ }\sC(s.\cA)
\xrightarrow{\ \sC(\sU_\cA)\ }\sC(\sC(\cA)).
$$
In the same vein, theorem \ref{th_degen-plus-norm} yields
natural decompositions of additive functors :
\set\begin{equation}\label{eq_tot-this}
\sU^2_\cA=(\sC(\sN_\cA)\oplus\sC(\sD_\cA))\circ
(\sN_{s.\cA}\oplus\sD_{s.\cA})
\end{equation}
and if we let
$$
\sN^2_\cA:=\sC(\sN_\cA)\circ\sN_{s.\cA}:s^2.\cA\to\sC(\sC(\cA))
\qquad
A\mapsto N_{\bullet\bullet}A
$$
the natural morphism
$$
\si^A_\bullet:\Tot\,N_{\bullet\bullet}A\to\Tot\,A_{\bullet\bullet}
$$
is a homotopy equivalence.

\sset\subsubsection{}\label{subsec_AW-and-Sh}
We wish next to exhibit two natural transformations
$$
A^\Delta_\bullet\xrightarrow{\ \AW_\bullet^A\ }
\Tot A_{\bullet\bullet}\xrightarrow{\ \Sh_\bullet^A\ }
A^\Delta_\bullet
\qquad
\text{for every $A\in\Ob(s^2.\cA)$}.
$$
Namely, for every $n\in\N$, define :
\begin{itemize}
\item
$\Sh_n^A$ as the sum, for all $p,q\in\N$ such that $p+q=n$,
of the {\em shuffle maps}
$$
\Sh^A_{p,q}:=\sum_{\mu\nu}\eps_{\mu\nu}\cdot
A[\eta_{\nu_1}\circ\cdots\circ\eta_{\nu_q},
\eta_{\mu_1}\circ\cdots\circ\eta_{\mu_p}]:
A_{p,q}\to A_{n,n}
$$
where the sum ranges over the {\em shuffle permutations}
$(\mu,\nu)$ of type $(p,q)$ of the set $\{0,\dots,p+q-1\}$
(these are the permutations described in \cite[\S4.3.15]{Ga-Ra}),
and $\eps_{\mu\nu}$ is the sign of the permutation $(\mu,\nu)$
\item
$\AW_n^A$ as the sum, for all $p,q\in\N$ such that $p+q=n$, of
the {\em Alexander-Whitney maps}
$$
\AW^A_{p,q}:=A[\eps^{q\vee}_{p,0},\eps^p_{q,0}]:A_{n,n}\to A_{p,q}
$$
\end{itemize}
(notation of example \ref{ex_face-and-deg}(ii)). Clearly, for
every $p,q\in\N$ the rule $A\mapsto A_{p,q}$ defines a functor
$$
\bullet[p,q]:s^2.\cA\to\cA
$$
and the maps $\Sh^A_{p,q}$ and $\AW^A_{p,q}$ yield natural
transformations
$$
\Sh_{p,q}:\bullet[p,q]\Rightarrow\bullet[p+q,p+q]
\qquad
\AW_{p,q}:\bullet[p+q,p+q]\Rightarrow\bullet[p,q].
$$

\begin{proposition}\label{prop_shuffle}
With the notation of \eqref{subsec_shuffle}, the sequence
$(\Sh_n^A~|~n\in\N)$ defines a morphism of chain complexes
$$
\Sh_\bullet^A:
\Tot A_{\bullet\bullet}\to A^\Delta_\bullet
\qquad
\text{for every $A\in\Ob(s^2.\cA)$}.
$$
\end{proposition}
\begin{proof} For any bisimplicial object $A$ of $\cA$,
set $A_{-1,q}=A_{p,-1}:=0$ for every $p,q\in\Z$, and let
$$
A[\eps_0,\one_{[q]}]:A_{0,q}\to A_{-1,q}
\qquad
A[\one_{[p]},\eps_0]:A_{p,0}\to A_{p,-1}
\qquad
A[\eps_0,\eps_0]:A_{0,0}\to A_{-1,-1}
$$
be the zero maps; likewise, let
$\Sh^A_{-1,q}:A_{-1,q}\to A_{q-1,q-1}$ and
$\Sh^A_{p,-1}:A_{p,-1}\to A_{p-1,p-1}$ be the zero maps,
and notice that also $\Sh^A_{0,0}$ is the zero map.
We define 
$$
d^{A,h}_{p,q}:=
\sum_{i=0}^p(-1)^i\cdot A[\eps_i,\one_{[q]}]:A_{p,q}\to A_{p-1,q}
\qquad
d^{A,v}_{p,q}:=\sum_{i=0}^q(-1)^j\cdot
A[\one_{[p]},\eps_i]:A_{p,q}\to A_{p,q-1}
$$
so, for every $n\in\N$ the differential $d_n$ in degree $n$ of
$\Tot A_{\bullet\bullet}$ is the sum of the maps
$$
d^{A,h}_{p,q}+(-1)^p\cdot d^{A,v}_{p,q}:A_{p,q}\to A_{p-1,q}\oplus A_{p,q-1}
\qquad
\text{for all $p,q\in\N$ such that $p+q=n$}
$$
whereas the differential of $A^\Delta_\bullet$ is the morphism
$$
d^A_n:=\sum^n_{i=1}A[\eps_i,\eps_i]:A_{n,n}\to A_{n-1,n-1}
$$
and we have to check the identity
\set\begin{equation}\label{eq_a-chain-map}
d^A_{p+q}\circ\Sh^A_{p,q}=
\Sh^A_{p-1,q}\circ d^{A,h}_{p,q}+(-1)^p\cdot\Sh^A_{p,q-1}\circ d^{A,v}_{p,q}
\qquad
\text{for all $p,q\in\N$}.
\end{equation}
Now, we set $B:=\gamma_1A$, $C:=\gamma_2A$, $D:=\gamma_2B$
(notation of \eqref{subsec_bisimplex}), but for the purpose
of this proof, we shall modify the differentials of the
double complexes $B_{\bullet\bullet}$ and $C_{\bullet\bullet}$
in certain low degrees : namely, we define
$$
d^{B,h}_{0,q}:=A[\eps_0,\one_{[q]}]
\qquad
d^{C,v}_{p,0}:=A[\one_{[p]},\eps_0]
\qquad
\text{for every $p,q\in\N$}.
$$

\begin{claim}\label{cl_crazy-combinatorics}
With the foregoing notation, the following holds :
\begin{enumerate}
\item
$\Sh^A_{p,q}=
(-1)^q\cdot \Sh^D_{p-1,q}\circ A[\one_{[p]},\eta_q]+
\Sh^D_{p,q-1}\circ A[\eta_p,\one_{[q]}]$ \ \ \
for every $p,q\in\N$.
\item
$A[\eps_{p+q-1},\eps_{p+q-1}]\circ\Sh^D_{p-1,q-1}=
\Sh^A_{p-1,q-1}\circ A[\eps_p,\eps_q]$ \ \ \ for every $p,q>0$.
\item
$A[\eps_{p+q},\eps_{p+q}]\circ\Sh^A_{p,q}=
(-1)^q\cdot\Sh^A_{p-1,q}\circ A[\eps_p,\one_{[q]}]+
\Sh^A_{p,q-1}\circ A[\one_{[p]},\eps_q]$ \ \ \ for every $p,q\in\N$.
\end{enumerate}
\end{claim}
\begin{pfclaim}(i): First, we notice the identities :
\set\begin{equation}\label{eq_with-aligned-for-A}
{\begin{aligned}
A[\one_{[p+q]},\eta_{p+q}]\circ\Sh^A_{p,q}=\, &
\Sh^C_{p,q}\circ A[\one_{[p]},\eta_q] \\
A[\eta_{p+q},\one_{[p+q]}]\circ\Sh^A_{p,q}=\, &
\Sh^B_{p,q}\circ A[\eta_p,\one_{[q]}].
\end{aligned}}
\qquad
\text{for every $p,q\in\N$}
\end{equation}
that are deduced from the identities
$$
\begin{aligned}
\eta_{\mu_1}\circ\cdots\circ\eta_{\mu_p}\circ\eta_{p+q}=\, &
\eta_q\circ\eta_{\mu_1}\circ\cdots\circ\eta_{\mu_p} \\
\eta_{\nu_1}\circ\cdots\circ\eta_{\nu_q}\circ\eta_{p+q}=\, &
\eta_p\circ\eta_{\nu_1}\circ\cdots\circ\eta_{\nu_q} 
\end{aligned}
$$
which in turn follow from the simplicial identities for
degeneracy maps. By applying the first (resp. second)
identity \eqref{eq_with-aligned-for-A} with $A$ replaced
by $B$ (resp. by $C$) and $p$ replaced by $p-1$ (resp.
with $q$ replaced by $q-1$) we get :
\set\begin{equation}\label{eq_with-aligned}
{\begin{aligned}
A[\one_{[p+q]},\eta_{p+q-1}]\circ\Sh^B_{p-1,q}=\, &
\Sh^D_{p-1,q}\circ A[\one_{[p]},\eta_q] \\
A[\eta_{p+q-1},\one_{[p+q]}]\circ\Sh^C_{p,q-1}=\, &
\Sh^D_{p,q-1}\circ A[\eta_p,\one_{[q]}].
\end{aligned}}
\qquad
\text{for every $p,q\in\N$}.
\end{equation}
Now, suppose first that $p,q>0$, and let $(\mu,\nu)$ be any
$(p,q)$-shuffle of $\{0,\dots,p+q-1\}$; then either $\mu_p=p+q-1$
or $\nu_q=p+q-1$. In the first (resp. second) case, after
removing $\mu_p$ (resp. $\nu_q$) we get a $(p-1,q)$-shuffle
(resp. a $(p,q-1)$-shuffle) $(\bar\mu,\bar\nu)$ with
$$
\eps_{\bar\mu\bar\nu}=(-1)^q\cdot\eps_{\mu\nu}
\qquad
\text{(resp. $\eps_{\bar\mu\bar\nu}=\eps_{\mu\nu}$)}.
$$
However, the first (resp. second) left-hand side of
\eqref{eq_with-aligned} contains precisely all the
terms of the first (resp. second) type occurring in the
definition of $\Sh^A_{p,q}$, so we get (i) in this case.
The cases where either $p=0$ or $q=0$ can be dealt with
by a similar, but simpler, argument.

(ii) follows by naturality of $\Sh_{p-1,q-1}$, applied
to the morphism $g^{(1)}_A\circ g^{(2)}_B:D\to A$ from
\eqref{subsec_bisimplex}.

(iii): The case where $p=q=0$ is obvious, and the other cases
follow by composing both sides of (i) with $A[\eps_{p+q},\eps_{p+q}]$,
applying (ii), and recalling that
$\eta_q\circ\eps_q=\one_{[q]}$ : details left to the reader. 
\end{pfclaim}

Now, a simple inspection shows that \eqref{eq_a-chain-map}
follows from claim \ref{cl_crazy-combinatorics}(iii) and the
following 

\begin{claim} $d^D_{p+q-1}\circ\Sh^A_{p,q}=
\Sh^A_{p-1,q}\circ d^{B,h}_{p-1,q}+
(-1)^p\cdot\Sh^A_{p,q-1}\circ d^{C,v}_{p,q-1}$ \ \ for every
$p,q\in\N$.
\end{claim}
\begin{pfclaim}[] Consider first the case where $p=0$.
If $q\leq 1$, it is easily seen that both sides of the
stated identity vanish; if $q\geq 2$, the left-hand side is
$$
d^D_{q-1}\circ A[\eta_0\circ\cdots\circ\eta_{q-1},\one_{[q]}]=
\sum_{i=0}^{q-1}(-1)^i\cdot
A[\eta_0\circ\cdots\circ\eta_{q-1}\circ\eps_i,\eps_i]
$$
and the right-hand side is
$$
A[\eta_0\circ\cdots\circ\eta_{q-2},\one_{[q]}]\circ d^{C,v}_{0,q-1}=
\sum_{i=0}^{q-1}(-1)^i\cdot
A[\eta_0\circ\cdots\circ\eta_{q-2},\eps_i]
$$
so in this case the assertion comes down to the obvious
identity
$$
\eta_0\circ\cdots\circ\eta_{q-1}\circ\eps_i=
\eta_0\circ\cdots\circ\eta_{q-2}:[q-2]\to[0]
\qquad
\text{for every $i=0,\dots,q-1$}.
$$
Likewise we deal with the case where $q=0$. It follows
already that \eqref{eq_a-chain-map} holds for these values
of $(p,q)$, and for every bisimplicial object $A$ of $\cA$;
especially, we get
\set\begin{equation}\label{eq_p-is-zero}
d^D_q\circ\Sh^D_{0,q}=\Sh^D_{0,q-1}\circ d^{D,v}_{0,q}
\qquad
\text{for every $q\in\N$}.
\end{equation}
Next, for $p=q=1$, a direct computation shows that both
sides equal
$$
A[\eps_0\circ\eta_0,\one_{[1]}]-A[\one_{[1]},\eps_0\circ\eta_0]:
A_{1,1}\to A_{1,1}.
$$
We prove now, by induction on $q$, that the assertion
holds for $p=1$. This is already known for $q\leq 1$,
so suppose that $r>1$, and that the assertion holds
for $p=1$ and every $q<r$. The latter implies that
also \eqref{eq_a-chain-map} holds for these values of
$(p,q)$, and for every bisimplicial object $A$ of $\cA$;
especially, we get
\set\begin{equation}\label{eq_p-is-one}
d^D_r\circ\Sh^D_{1,r-1}=\Sh^D_{0,r-1}\circ d^{D,h}_{1,r-1}-
\Sh^D_{1,r-2}\circ d^{D,v}_{1,r-1}.
\end{equation}
On the other hand, claim \ref{cl_crazy-combinatorics}(i)
says that
$$
d^D_r\circ\Sh^A_{1,r}=
(-1)^r\cdot d^D_r\circ\Sh^D_{0,r}\circ A[\one_{[1]},\eta_r]+
d^D_r\circ\Sh^D_{1,r-1}\circ A[\eta_1,\one_{[r]}]
$$
Combining with \eqref{eq_p-is-zero} and \eqref{eq_p-is-one}
we obtain
$$
d^D_r\circ\Sh^A_{1,r}=(-1)^r\cdot
\Sh^D_{0,r-1}\circ d^{D,v}_{0,r}\circ A[\one_{[1]},\eta_r]+
(\Sh^D_{0,r-1}\circ d^{D,h}_{1,r-1}-\Sh^D_{1,r-2}\circ d^{D,v}_{1,r-1})
\circ A[\eta_1,\one_{[r]}].
$$
However, we have
$$
d^{D,h}_{1,r-1}\circ A[\eta_1,\one_{[r]}]=
A[\eta_1\circ\eps_0,\one_{[r]}]-A[\one_{[1]},\one_{[r]}]=
A[\eps_0\circ\eta_0,\one_{[r]}]-\one_{A[1,r]}
$$
so we can rewrite $d^D_r\circ\Sh^A_{1,r}=\phi_1+\phi_2-\phi_3$,
where
$$
\begin{aligned}
\phi_1:=\, & \Sh^D_{0,r-1}\circ
((-1)^r\cdot d^{D,v}_{0,r}\circ A[\one_{[1]},\eta_r]-\one_{A[1,r]})=
(-1)^r\cdot\Sh^D_{0,r-1}\circ A[\one_{[1]},\eta_{r-1}]\circ d^{C,v}_{1,r-1} \\
\phi_2:=\, & \Sh^D_{0,r-1}\circ A[\eps_0\circ\eta_0,\one_{[r]}] \\
\phi_3:=\, & \Sh^D_{1,r-2}\circ d^{D,v}_{1,r-1}\circ A[\eta_1,\one_{[r]}].
\end{aligned}
$$
On the other hand, using claim \ref{cl_crazy-combinatorics}(i),
we can compute
$$
\begin{aligned}
\Sh^A_{0,r}\circ d^{B,h}_{0,r-1}=\, &
\Sh^D_{0,r-1}\circ A[\eta_0,\one_{[r]}]\circ d^{B,h}_{0,r-1}=\phi_2 \\
\Sh^A_{1,r-1}\circ d^{C,v}_{1,r-1}=\, &
((-1)^{r-1}\cdot \Sh^D_{0,r-1}\circ A[\one_{[1]},\eta_{r-1}]+
\Sh^D_{1,r-2}\circ A[\eta_1,\one_{[r-1]}])\circ d^{C,v}_{1,r-1} \\
=\, & \phi_3-\phi_1
\end{aligned}
$$
which shows that the claim holds for $p=1$ and $q=r$, and
concludes the induction.

Lastly, we prove the claim for every $p,q\in\N$, by induction
on $p+q$; notice that the foregoing already shows that the
assertion holds whenever $p+q\leq 2$. Thus, let $r>2$, and
suppose that the claim is already known for every pair
$(p,q)$ such that $p+q<r$; then also \eqref{eq_a-chain-map}
holds for such values of $p$ and $q$, and especially we get
\set\begin{equation}\label{eq_combine-crazy}
d^D_{p+q}\circ\Sh^D_{p,q}=
\Sh^D_{p-1,q}\circ d^{D,h}_{p,q}+(-1)^p\cdot\Sh^D_{p,q-1}\circ d^{D,v}_{p,q}
\qquad
\text{whenever $p+q<r$}.
\end{equation}
Let $(p',q')$ be a pair such that $p'+q'=r$; combining
\eqref{eq_combine-crazy} with claim
\ref{cl_crazy-combinatorics}(i), we get
$$
\begin{aligned}
d^D_{p'+q'-1}\circ\Sh^A_{p',q'}=\, & (-1)^{q'}\cdot
(\Sh^D_{p'-2,q'}\circ d^{D,h}_{p'-1,q'}-(-1)^{p'}\cdot
\Sh^D_{p'-1,q'-1}\circ d^{D,v}_{p'-1,q'})\circ A[\one_{[p']},\eta_{q'}] \\
& +(\Sh^D_{p'-1,q'-1}\circ d^{D,h}_{p',q'-1}+(-1)^{p'}\cdot
\Sh^D_{p',q'-2}\circ d^{D,v}_{p',q'-1})
\circ A[\eta_{p'},\one_{[q']}].
\end{aligned}
$$
On the other hand, after noticing that
$$
\begin{aligned}
d^{D,h}_{p',q'-1}\circ A[\eta_{p'},\one_{[q']}]-
(-1)^{p'}\cdot\one_{A[p',q']}=\, &
A[\eta_{p'-1},\one_{[q']}]\circ d^{B,h}_{p'-1,q'} \\
d^{D,v}_{p'-1,q'}\circ A[\one_{[p']},\eta_{q'}]-
(-1)^{q'}\cdot\one_{A[p',q']}=\, &
A[\one_{[p']},\eta_{q'-1}]\circ d^{C,v}_{p',q'-1} 
\end{aligned}
$$
we may apply claim \ref{cl_crazy-combinatorics}(i), to
compute
$$
\begin{aligned}
\Sh^A_{p'-1,q'}\circ d^{B,h}_{p'-1,q'}=\, &
((-1)^{q'}\cdot \Sh^D_{p'-2,q'}\circ A[\one_{[p'-1]},\eta_{q'}]+
\Sh^D_{p'-1,q'-1}\circ A[\eta_{p'-1},\one_{[q']}])\circ d^{B,h}_{p'-1,q'} \\
=\, & (-1)^{q'}\cdot \Sh^D_{p'-2,q'}\circ d^{D,h}_{p'-1,q'}\circ
A[\one_{[p']},\eta_{q'}]-(-1)^{p'}\cdot\Sh^D_{p'-1,q'-1} \\
& +\Sh^D_{p'-1,q'-1}\circ d^{D,h}_{p',q'-1}\circ A[\eta_{p'},\one_{[q']}]
\end{aligned}
$$
and
$$
\begin{aligned}
\Sh^A_{p',q'-1}\circ d^{C,v}_{p',q'-1}=\, &
((-1)^{q'-1}\cdot \Sh^D_{p'-1,q'-1}\circ A[\one_{[p']},\eta_{q'-1}] \\
\, & +
\Sh^D_{p',q'-2}\circ A[\eta_{p'},\one_{[q'-1]}])\circ d^{C,v}_{p',q'-1} \\
=\,& (-1)^{q'-1}\cdot\Sh^D_{p'-1,q'-1}\circ d^{D,v}_{p'-1,q'}
\circ A[\one_{[p']},\eta_{q'}]+\Sh^D_{p'-1,q'-1} \\
\, & +\Sh^D_{p',q'-2}\circ d^{D,v}_{p',q'-1}\circ A[\eta_{p'},\one_{[q']}].
\end{aligned}
$$
Finally, the sought identity for the pair $(p',q')$ follows
by comparing the last two identities with the previous one for
$d^D_{p'+q'-1}\circ\Sh^A_{p',q'}$, and this concludes the proof
of the inductive step.
\end{pfclaim}
\end{proof}

\begin{proposition}\label{prop_Alex-Whitney}
With the notation of \eqref{subsec_shuffle}, the sequence
$(\AW_n^A~|~n\in\N)$ defines a morphism of chain complexes
$$
\AW_\bullet^A:
A^\Delta_\bullet\to\Tot A_{\bullet\bullet}
\qquad
\text{for every $A\in\Ob(s^2.\cA)$}.
$$
\end{proposition}
\begin{proof} Denote by $d^T_\bullet$ the differential of
$\Tot A_{\bullet\bullet}$, and keep the notation of the proof
of proposition \ref{prop_shuffle}; we have to check the identity
\set\begin{equation}\label{eq_proj-hand-side}
\AW_{n-1}\circ d^A_n=d^T_n\circ\AW_n
\qquad
\text{for every $n\in\N$}.
\end{equation}
However, say that $p,q\in\N$ and $p+q=n$; the projection
of the right-hand side of \eqref{eq_proj-hand-side} on the
direct summand $A_{p-1,q}$ of $(\Tot A_{\bullet\bullet})_{n-1}$
equals
\set\begin{equation}\label{eq_unwind-this}
d^{A,h}_{p,q}\circ\AW_{p,q}+
(-1)^{p-1}\cdot d^{A,v}_{p-1,q+1}\circ\AW_{p-1,q+1}.
\end{equation}
Unwinding the definition, \eqref{eq_unwind-this} is found to be
$$
\sum_{i=0}^p(-1)^i\cdot
A[\eps^{q\vee}_{p,0}\circ\eps_i,\eps^p_{q,0}]-
(-1)^p\cdot\sum_{i=0}^{q+1}(-1)^i\cdot
A[\eps^{q+1\vee}_{p-1,0},\eps^{p-1}_{q+1,0}\circ\eps_i]
$$
which we may rewrite as
\set\begin{equation}\label{eq_ci-sta}
\sum_{i=0}^{p-1}(-1)^i\cdot
A[\eps^{q\vee}_{p,0}\circ\eps_i,\eps^p_{q,0}]+
\sum_{i=p}^{p+q}(-1)^i\cdot
A[\eps^{q+1\vee}_{p-1,0},\eps^{p-1}_{q+1,0}\circ\eps_{i-p+1}].
\end{equation}
On the other hand, the projection of the left-hand side of
\eqref{eq_proj-hand-side} onto $A_{p-1,q}$ equals
$$
\sum_{i=0}^n(-1)^i\cdot
A[\eps_i\circ\eps^{q\vee}_{p-1,0},\eps_i\circ\eps^{p-1}_{q,0}].
$$
To compare the latter with \eqref{eq_ci-sta}, it suffices
to remark that
$$
\eps_i\circ\eps^{q\vee}_{p-1,0}=\left\{
\begin{array}{ll}
\eps^{q\vee}_{p,0}\circ\eps_i & \qquad \text{if $i<p$} \\
\eps^{q+1\vee}_{p-1,0} & \qquad \text{if $i\geq p$}
\end{array}\right.
\qquad
\text{and}
\qquad
\eps_i\circ\eps_{q,0}^{p-1}=\left\{
\begin{array}{ll}
\eps^p_{q,0} & \qquad \text{if $i<p$} \\
\eps^{p-1}_{q+1,0}\circ\eps_i & \qquad \text{if $i\geq p$}
\end{array}\right.
$$
which are all deduced from the simplicial identities for the
$\eps_i$. The proposition follows.
\end{proof}

\sset\subsubsection{}\label{subsec_harlem-shuffle}
Propositions \ref{prop_shuffle} and \ref{prop_Alex-Whitney}
yield natural transformations of functors
$$
\xymatrix{
\Tot\,\sU^2_\cA \ar@<.5ex>[rr]^-{\ \Sh_\bullet\ }
& & \ar@<.5ex>[ll]^-{\ \AW_\bullet\ } \sU_\cA\circ\sDelta_\cA.
}$$
Now, denote by $\shh_\cA:\sC(\cA)\to\Hot(\cA)$ the natural
functor (this induces the identity on the objects, and
the projection on the group of morphisms); there follow
natural transformations
$$
\xymatrix{
\shh_\cA\circ\Tot\,\sU^2_\cA
\ar@<.5ex>[rr]^-{\ \shh_\cA*\Sh_\bullet\ }
& & \ar@<.5ex>[ll]^-{\ \shh_\cA*\AW_\bullet\ }
\shh_\cA\circ\sU_\cA\circ\sDelta_\cA.
}
\qquad
\text{between functors $s^2.\cA\to\Hot(\cA)$.}
$$

\begin{theorem}[Eilenberg-Zilber-Cartier]
\label{th_Eilenberg-Zilber}
With the notation of \eqref{subsec_harlem-shuffle}, we have :
\begin{enumerate}
\item
$\shh_\cA*\AW_\bullet$ and $\shh_\cA*\Sh_\bullet$ are mutually
inverse isomorphisms of functors.
\item
More precisely, there exist natural modifications (see definition
{\em\ref{def_modification}})
$$
\AW_\bullet\circ\Sh_\bullet\leadsto\one_{\Tot\,\sU^2_\cA}
\qquad
\Sh_\bullet\circ\AW_\bullet\leadsto\one_{\sU_\cA\circ\Delta_\cA}.
$$
\end{enumerate}
\end{theorem}
\begin{proof} Let
$\sp^A_\bullet:\Tot\,A_{\bullet\bullet}\to\Tot\,N_{\bullet\bullet}A$
be the projection whose kernel is the sum of the remaining
three direct summands in the decomposition
$\Tot\,\eqref{eq_tot-this}$; explicitly, in each degree $n\in\N$,
this kernel is the sum of the subobjects
$\Img\,A[\eta_i,\one_{[q]}]$ and $\Img\,A[\one_{[p]},\eta_j]$
of $A_{p,q}$, for every $i=0,\dots,p-1$, $j=0,\dots,q-1$, and
every $p,q\in\N$ such that $p+q=n$. Likewise, let
$\sj^A_\bullet:N_\bullet A^\Delta\to A_\bullet^\Delta$ and
$\sq^A_\bullet:A_\bullet^\Delta\to N_\bullet A^\Delta$ be as
in \eqref{subsec_sj-and-sq}; theorem \ref{th_degen-plus-norm}
and the discussion of \eqref{subsec_iterate-U} show that
the pairs $(\sp^A_\bullet,\si^A_\bullet)$ and $(\sj_\bullet^A,\sq^A_\bullet)$
induce mutually inverse isomorphisms in $\Hot(\cA)$, so it suffices
to check that the same holds for the compositions
$$
\bar\Sh{}^A_\bullet:=\sq_\bullet^A\circ\Sh^A_\bullet\circ\si_\bullet^A:
\Tot\,N_{\bullet\bullet}A\to N_\bullet A^\Delta
\qquad
\bar\AW{}^A_\bullet:=\sp_\bullet^A\circ\AW^A_\bullet\circ\sj_\bullet^A:
N_\bullet A^\Delta\to\Tot\,N_{\bullet\bullet}A
$$
for every bisimplicial object $A$ of $\cA$. However, we have

\begin{claim}\label{cl_first-modif}
$\sp^A_\bullet\circ\AW^A_\bullet\circ\Sh^A_\bullet\circ\si^A_\bullet
=\one_{\Tot\,N_{\bullet\bullet}A}$.
\end{claim}
\begin{pfclaim} More precisely, say that $p+q=n$, consider
the composition
$$
f_{p,q}:A_{p,q}\xrightarrow{\ \Sh^A_{p,q}\ }A_{p+q,p+q}
\xrightarrow{\ \AW^A_n\ }(\Tot\,A_{\bullet\bullet})_n
$$
and let $g_{p,q}:A_{p,q}\to(\Tot\,A_{\bullet\bullet})_n$ be
the inclusion map; we shall show that
$f_{p,q}-g_{p,q}\subset\Ker\,\sp^A_n$. Indeed, we have
$$
f_{p,q}=\sum_{i=0}^n\sum_{(\mu,\nu)}\eps_{\mu\nu}\cdot
A[\eta_{\nu_1}\circ\cdots\circ\eta_{\nu_q}\circ\eps^{n-i\vee}_{i,0},
\eta_{\mu_1}\circ\cdots\circ\eta_{\mu_p}\circ\eps^i_{n-i,0}]
$$
(notation of \eqref{subsec_AW-and-Sh}, and $\eps^i_{n-i,0}$
is defined as in the proof of proposition \ref{prop_Alex-Whitney}).
However, if $i<p$, then
$\eta_{\mu_1}\circ\cdots\circ\eta_{\mu_p}\circ\eps^i_{n-i,0}$
is of the form $\tau\circ\eta_{\mu_p}$ for some
$\tau:[p+q-i-1]\to[q]$, so the corresponding term does lie in
$\Ker\,\sp^A_n$. If $i>p$, then
$\eta_{\nu_1}\circ\cdots\circ\eta_{\nu_q}\circ\eps^{n-i\vee}_{i,0}$
is of the form $\tau\circ\eta_k$ for some $k\leq\nu_q$ and
some $\tau:[i-1]\to[p]$, so this term likewise lies in
$\Ker\,\sp^A_n$. It remains to consider the terms with $i=p$;
however, a simple inspection shows that
$\eta_{\nu_1}\circ\cdots\circ\eta_{\nu_q}\circ\eps^{n-p\vee}_{p,0}$
is of the same form as above, unless $\nu_1$ is either $p-1$ or $p$
(details left to the reader); furthermore, if $\nu_1=p-1$,
then $\mu_p\geq p$, in which case
$\eta_{\mu_1}\circ\cdots\circ\eta_{\mu_p}\circ\eps^i_{n-i,0}$
is of the form described above. In all these cases, the
corresponding term again lies in $\Ker\,\sp^A_n$. So, it
remains only to consider the single case where $(\mu,\nu)$
is the identity permutation, whose sign equals $1$; in
this case, the corresponding term is none else than
$A[\one_{[p]},\one_{[q]}]$, whence the claim.
\end{pfclaim}

Claim \ref{cl_first-modif} and theorem
\ref{th_degen-plus-norm}(ii) already yield the existence
of the sought modification
$\AW_\bullet\circ\Sh_\bullet\leadsto\one_{\Tot\,\sU^2_\cA}$.
Next we define, for every $n\in\N$, a natural transformation
$$
s_n:\bullet[n]\circ\sDelta_\cA\Rightarrow
\bullet[n+1]\circ\sDelta_\cA
$$
(notation of \eqref{subsec_simplicial-object}; so $s_n^A$ is
a morphism $A[n,n]\to A[n+1,n+1]$ for every object
$A$ of $s^2.\cA$). The construction is by induction on $n$ : for
$n=0$ we let $s_0^A:A[0,0]\to A[1,1]$ be the zero morphism;
for $n>0$ we set
$$
s^A_n:=\Sh^{A'}_n\circ\AW^{A'}_n\circ A[\eta_0,\eta_0]-s_{n-1}^{A'}
\qquad
\text{where $A':=(\gamma_2\circ\gamma_1(A^\vee))^\vee$}
$$
(notation of \eqref{subsec_bisimplex} and \eqref{subsec_front-back} :
explicitly, we have $A'[p,q]:=A[p+1,q+1]$ for every $p,q\in\N$,
and $A'[\eps_i,\eps_j]:=A[\eps_{i+1},\eps_{j+1}]$, and likewise
for $A'[\eps_i,\eta_j]$, $A'[\eta_i,\eps_j]$ and $A'[\eta_i,\eta_j]$,
for every face and degeneracy map of $\Delta$).
We have

\begin{claim}\label{cl_second-modif}
The system $(\sq^A_{n+1}\circ s^A_n\circ\sj^A_n~|~n\in\N)$ is a
homotopy
$\sq^A_\bullet\circ\Sh^A_\bullet\circ\AW^A_\bullet\circ\sj^A_\bullet
\Rightarrow\one_{N_\bullet A^\Delta}$.
\end{claim}
\begin{pfclaim} Denote by $d^A_\bullet$ the differential of
$A^\Delta_\bullet$, and let $d^A_0:A^\Delta_0\to 0$ and
$s^A_{-1}:0\to A^\Delta_0$ be the zero morphisms.  We check,
more precisely, that
$$
\sq^A_n\circ(d^A_{n+1}\circ s^A_n+s^A_{n-1}\circ d^A_n)=
\sq^A_n\circ\Sh^A_n\circ\AW^A_n-\sq^A_n
\qquad
\text{for every $n\in\N$}.
$$
We argue by induction on $n$, and the assertion is clear for
$n=0$. For $n=1$, notice that
$$
\Sh^A_1\circ\AW^A_1=
A[\eps_1\circ\eta_0,\one_{[1]}]+A[\one_{[1]},\eps_0\circ\eta_0]:
A_{1,1}\to A_{1,1}.
$$
It follows that
$$
s^A_1=\Sh^{A'}_1\circ\AW^{A'}_1\circ A[\eta_0,\eta_0]
=(A[\eps_1\circ\eta_0\circ\eta_0,\eta_0]+A[\eta_0,\eta_1])
$$
whence
$$
d^A_2\circ s^A_1=\Sh^A_1\circ\AW^A_1-\one_{A[1,1]}+
A[\eps_1\circ\eps_1\circ\eta_0,\eps_1\circ\eta_0]
$$
(details left to the reader). Since
$\Img\,A[\eta_0,\eta_0]\subset\Ker\,\sq^A_2$, the assertion
follows in this case. Next, suppose that $r>1$, and that the
sought identity is already known for every $n<r$, and every
bisimplicial object $A$; especially, it holds for $A'$, so
if we let $d^{A'}_\bullet$ be the differential of $A'_\bullet$,
we get
\set\begin{equation}\label{eq_yoga}
\sq^{A'}_{r-1}\circ(
d^{A'}_r\circ s^{A'}_{r-1}+s^{A'}_{r-2}\circ d^{A'}_{r-1})=
\sq^{A'}_{r-1}\circ\Sh^{A'}_{r-1}\circ
\AW^{A'}_{r-1}-\sq^{A'}_{r-1}.
\end{equation}
On the other hand, after noticing that
$$
d^{A'}_r\circ A[\eta_0,\eta_0]=
\one_{A[r,r]}-A[\eta_0,\eta_0]\circ d^{A'}_{r-1}
$$
we may compute
$$
\begin{aligned}
d^{A'}_r\circ s^A_r=\, &
d^{A'}_r\circ\Sh^{A'}_r\circ\AW^{A'}_r\circ A[\eta_0,\eta_0]-
d^{A'}_r\circ s^{A'}_{r-1} \\
=\, &
\Sh^{A'}_r\circ\AW^{A'}_r\circ d^{A'}_r\circ A[\eta_0,\eta_0]-
d^{A'}_r\circ s^{A'}_{r-1} \\
=\, &
\Sh^{A'}_r\circ\AW^{A'}_r-
\Sh^{A'}_r\circ\AW^{A'}_r\circ A[\eta_0,\eta_0]\circ d^{A'}_{r-1}-
d^{A'}_r\circ s^{A'}_{r-1}
\end{aligned}
$$
which, combined with \eqref{eq_yoga}, implies
\set\begin{equation}\label{eq_mardi-yoga}
\sq^{A'}_{r-1}\circ(d^{A'}_r\circ s^A_r+s^A_{r-1}\circ d^{A'}_{r-1})=
-\sq^{A'}_{r-1}.
\end{equation}
Furthermore, notice the natural morphism in $s^2.\cA$
$$
(g_{A^\vee}^{(1)}\circ g^{(2)}_{\gamma_1A^\vee})^\vee:A'\to A
$$
supplied by \eqref{subsec_bisimplex}; explicitly, for every
$p,q\in\N$, this morphism is given by the discarded face
operator $A[\eps_0,\eps_0]:A[p+1,q+1]\to A[p,q]$. Then,
the naturality of $s_\bullet$, $\Sh_\bullet$ and $\AW_\bullet$
implies
$$
\begin{aligned}
A[\eps_0,\eps_0]\circ s^A_r=\, &
A[\eps_0,\eps_0]\circ
\Sh^{A'}_r\circ\AW^{A'}_r\circ A[\eta_0,\eta_0]-
A[\eps_0,\eps_0]\circ s_{r-1}^{A'} \\
=\, &
\Sh^A_r\circ\AW^A_r\circ A[\eta_0\circ\eps_0,\eta_0\circ\eps_0]-
s_{r-1}^A\circ A[\eps_0,\eps_0] \\
=\, &
\Sh^A_r\circ\AW^A_r-s_{r-1}^A\circ A[\eps_0,\eps_0].
\end{aligned}
$$
Lastly, recalling that $d^A_r=A[\eps_0,\eps_0]-d^{A'}_{r-1}$,
the latter identity can be added to \eqref{eq_mardi-yoga}, to
deduce
$$
\sq^{A'}_{r-1}\circ(d_{r+1}^A\circ s^A_r+s^A_{r-1}\circ d^A_r)=
\sq^{A'}_{r-1}\circ\Sh^A_r\circ\AW^A_r-\sq^{A'}_{r-1}.
$$
Now, to prove the assertion in degree $r$, it suffices
to observe that $\sq^A_r$ factors through $\sq^{A'}_{r-1}$.
\end{pfclaim}

Claim \ref{cl_second-modif} and theorem \ref{th_degen-plus-norm}(ii)
supply the second sought modification, and conclude the proof
of the theorem.
\end{proof}

\sset\subsubsection{}\label{subsec_shuffle-commute}
Let $(\cA,\otimes)$ be an abelian tensor category, $A[\bullet]$
and $B[\bullet]$ two objects of $s.\cA$, and define the
bisimplicial objects $A\boxtimes B$ and $B\boxtimes A$, as well
as the simplicial objects $A\otimes B$ and $B\otimes A$  of $\cA$
as in remark \ref{rem_mixed-simpl-tensors}(iii). Notice that the
system of commutativity constraints $(\Psi_{A[n],B[n]}~|~n\in\N)$
amounts to an isomorphism
$$
\Psi_{A\otimes B}:A\otimes B\isom B\otimes A
\qquad
\text{in $s.\cA$}
$$
whence an isomorphism
$\Psi_{(A\otimes B)_\bullet}:(A\otimes B)_\bullet\isom(B\otimes A)_\bullet$
on the respective unnormalized complexes.

\begin{proposition}\label{prop_commu-shuffle}
With the notation of \eqref{subsec_shuffle-commute}, the
diagram of chain complexes
$$
\xymatrix{
\Tot(A\boxtimes B)_{\bullet\bullet}
\ar[rr]^-{\Psi^\bullet_{A_\bullet,B_\bullet}}
\ar[d]_{\Sh^{A\boxtimes B}_\bullet} & & \Tot(B\boxtimes A)_{\bullet\bullet}
\ar[d]^-{\Sh^{B\boxtimes A}_\bullet} \\
(A\otimes B)_\bullet \ar[rr]^-{\Psi_{(A\otimes B)_\bullet}} & &
(B\otimes A)_\bullet
}$$
commutes, where $\Psi^\bullet_{A_\bullet,B_\bullet}$ is the
commutativity constraint for the unnormalized chain complexes
$A_\bullet$ and $B_\bullet$, as in \eqref{eq_commut-constraint}.
\end{proposition}
\begin{proof} The assertion boils down to the identity
\set\begin{equation}\label{eq_boils-down-shuufle}
(-1)^{pq}\cdot\Sh^{B\boxtimes A}_{q,p}\circ\Psi_{A[p],B[q]}=
\Psi_{A[n],B[n]}\circ\Sh^{A\boxtimes B}_{p,q}
\end{equation}
for every $p,q\in\N$ with $p+q=n$. For the latter, suppose
first that $p=0$, in which case
$$
\Sh^{A\boxtimes B}_{p,q}=
A[\eta_0\circ\cdots\circ\eta_{q-1}]\otimes\one_{B[q]}
\qquad\text{and}\qquad
\Sh^{B\boxtimes A}_{q,p}=
\one_{B[q]}\otimes A[\eta_0\circ\cdots\circ\eta_{q-1}]
$$
from which we derive \eqref{eq_boils-down-shuufle}, using
the naturality of $\Psi$ (details left to the reader). Likewise
we argue for the case where $q=0$. For the general case, we
proceed by induction on $n$. The cases $n=0,1$ have already
been dealt with, so suppose $r\geq 2$, and that the sought
identity is already known for every pair of integers whose sum
is $<n$, and every objects $A,B$ of $s.\cA$. By the foregoing,
we may also assume that both $p,q>0$, and then claim
\ref{cl_crazy-combinatorics}(i) implies that
$$
\Sh^{A\boxtimes B}_{p,q}=
(-1)^q\cdot\Sh^{\gamma A\boxtimes\gamma B}_{p-1,q}\circ
(\one_{A[p]}\otimes B[\eta_q])+
\Sh^{\gamma A\boxtimes\gamma B}_{p,q-1}\circ
(A[\eta_p]\otimes\one_{B[q]}).
$$
However, it follows easily from remark \ref{rem_fun-exact-add}(iii)
that $\Psi$ is an additive functor in both of its arguments;
combining with the inductive assumption, we deduce that
$$
\begin{aligned}
\Psi_{A[n],B[n]}\circ\Sh^{A\boxtimes B}_{p,q}=\, &
(-1)^{pq}\cdot\Sh^{\gamma B\boxtimes\gamma A}_{q,p-1}\circ
\Psi_{A[p],B[q+1]}\circ(\one_{A[p]}\otimes B[\eta_q]) \\
& +(-1)^{p(q-1)}\cdot\Sh^{\gamma B\boxtimes\gamma A}_{q-1,p}\circ
\Psi_{A[p+1],B[q]}\circ(A[\eta_p]\otimes\one_{B[q]}) \\
=\, & (-1)^{pq}\cdot\Sh^{\gamma B\boxtimes\gamma A}_{q,p-1}
\circ(B[\eta_q]\otimes\one_{A[p]})\circ\Psi_{A[p],B[q]} \\
& +(-1)^{p(q-1)}\cdot\Sh^{\gamma B\boxtimes\gamma A}_{q-1,p}\circ
(\one_{B[q]}\otimes A[\eta_p])\circ\Psi_{A[p],B[q]} \\
=\, & (-1)^{pq}\cdot\Sh^{B\boxtimes A}_{q,p}\circ\Psi_{A[p],B[q]}
\end{aligned}
$$
where the last equality follows again from claim
\ref{cl_crazy-combinatorics}(i) and the additivity of $\Psi$.
\end{proof}

\sset\subsubsection{}\label{subsec_asso-shuffle}
The shuffle map is also {\em associative}, in the following sense.
Let
$$
A:=(A[p,q,r]~|~p,q,r\in\N)
$$
be any triple simplicial object of the abelian category $\cA$,
and denote by $A^{(1,2)}$ (resp. $A^{(2,3)}$) the diagonal
bisimplicial object of $\cA$ extracted from $A$ by the rule
$$
A^{(1,2)}[p,q]:=A[p,p,q]
\qquad
\text{(\ resp. $A^{(2,3)}[p,q]:=A[p,q,q]$\ )}
\qquad
\text{for every $p,q\in\N$}.
$$
Let also $A_{\bullet\bullet\bullet}$ be the triple chain complex
associated to $A$, and $A_{\bullet\bullet\bullet}^\Delta$ the
diagonal chain complex extracted from $A_{\bullet\bullet\bullet}$.
Moreover, denote by $A'$ (resp. $A''$) the bisimplicial object of
$s.\cA$ given by the rule :
$$
[p,q]\mapsto A[p,q,\bullet]
\quad
\text{(\ resp. $[p,q]\mapsto A[\bullet,p,q]$\ )}
\qquad
\text{for every $p,q\in\N$}.
$$

\begin{proposition}\label{prop_asso-shuffle}
With the notation of \eqref{subsec_asso-shuffle}, we have
a commutative diagram in $\sC(\cA)$
$$
\xymatrix{ \Tot(A_{\bullet\bullet\bullet}) \ar[rr] \ar[d]
& & \Tot(A_{\bullet\bullet}^{(1,2)}) \ar[d]^{\Sh^{A^{(1,2)}}_\bullet}  \\
\Tot(A_{\bullet\bullet}^{(2,3)}) \ar[rr]^-{\Sh^{A^{(2,3)}}_\bullet} & &
A_{\bullet\bullet\bullet}^\Delta
}$$
whose top horizontal (resp. left vertical) arrow is
obtained as the composition of $\Sh^{A'}_\bullet$ (resp.
$\Sh^{A''}_\bullet$) with the functor
$\Tot\circ\sC(\sU_\cA):\sC(s.\cA)\to\sC(\sC(\cA))\to\sC(\cA)$. 
\end{proposition}
\begin{proof} For every $p,q,r\in\N$, denote by
$$
\begin{aligned}
\Sh^{A'}_{p,q}[r]:\, & A[p,q,r]\to A[p+q,p+q,r]=A^{(1,2)}[p+q,r] \\
\Sh^{A''}_{q,r}[p]:\, & A[p,q,r]\to A^{(2,3)}[p,q+r]
\end{aligned}
$$
respectively the $[r]$-component of $\Sh^{A'}_{p,q}$ and the
$[p]$-component of $\Sh^{A''}_{q,r}$. The assertion boils down
to the identity :
$$
\Sh^{A^{(1,2)}}_{p+q,r}\circ\Sh^{A'}_{p,q}[r]=
\Sh^{A^{(2,3)}}_{p,q+r}\circ\Sh^{A''}_{q,r}[p]
\qquad
\text{for every $p,q,r\in\N$}.
$$
To check the latter, set
$$
D':=\gamma_2\circ\gamma_1A'
\qquad
D^{(1,2)}:=\gamma_2\circ\gamma_1(A^{(1,2)})
$$
and define likewise $D''$ and $D^{(2,3)}$ (notation of
\eqref{subsec_bisimplex}). Also, let $A[p,q,r]:=0$ whenever
one of the indices $p,q,r$ is $<0$, and define $\Sh_{p,q}$
to be the zero map, when either $p$ or $q$ is strictly
negative (cp. the proof of proposition \ref{prop_shuffle}).
We argue by induction on $n:=p+q+r$; the case $n=0$ is
trivial, so suppose that $n>0$, and that the assertion is
already known for all indices whose sum is $<n$, and all
triple simplicial objects of $\cA$. Notice as well that the
assertion trivially holds as well if any of the indices
$p,q,r$ is strictly negative, since in this case both sides
are the zero map. Hence, we may assume that $p,q,r\in\N$; in
this case, by applying claim \ref{cl_crazy-combinatorics}(i),
first to $A^{(1,2)}$ and then to $A'$, we compute
$$
\begin{aligned}
\Sh^{A^{(1,2)}}_{p+q,r}\circ\Sh^{A'}_{p,q}[r]=\, &
(-1)^r\cdot\Sh^{D^{(1,2)}}_{p+q-1,r}\circ A[\one_{[p+q]},\one_{[p+q]},\eta_r]
\circ\Sh^{A'}_{p,q}[r] \\
& +\Sh^{D^{(1,2)}}_{p+q,r-1}\circ A[\eta_{p+q},\eta_{p+q},\one_{[r]}]
\circ\Sh^{A'}_{p,q}[r] \\
=\, & (-1)^r\cdot\Sh^{D^{(1,2)}}_{p+q-1,r}\circ\Sh^{A'}_{p,q}[r+1]
\circ A[\one_{[p]},\one_{[q]},\eta_r] \\
& +\Sh^{D^{(1,2)}}_{p+q,r-1}\circ\Sh^{D'}_{p,q}[r]\circ
A[\eta_p,\eta_q,\one_{[r]}]
& \quad\text{(by \eqref{eq_with-aligned-for-A})} \\
=\, & (-1)^{q+r}\cdot\Sh^{D^{(1,2)}}_{p+q-1,r}\circ
\Sh^{D'}_{p-1,q}[r+1]\circ A[\one_{[p]},\eta_q,\eta_r] \\
& +(-1)^r\cdot\Sh^{D^{(1,2)}}_{p+q-1,r}\circ
\Sh^{D'}_{p,q-1}[r+1]\circ A[\eta_p,\one_{[q]},\eta_r] \\
& +\Sh^{D^{(1,2)}}_{p+q,r-1}\circ\Sh^{D'}_{p,q}[r]\circ
A[\eta_p,\eta_q,\one_{[r]}] \\
=\, & (-1)^{q+r}\cdot\Sh^{D^{(2,3)}}_{p-1,q+r}\circ\Sh^{D''}_{q,r}[p]
\circ A[\one_{[p]},\eta_q,\eta_r] \\
& +(-1)^r\cdot\Sh^{D^{(2,3)}}_{p,q+r-1}\circ\Sh^{D''}_{q-1,r}[p+1]
\circ A[\eta_p,\one_{[q]},\eta_r] \\
& +\Sh^{D^{(2,3)}}_{p,q+r-1}\circ\Sh^{D''}_{q,r-1}[p+1]
\circ A[\eta_p,\eta_q,\one_{[r]}]
\end{aligned}
$$
where the last identity holds by inductive assumption. On the
other hand, by applying claim \ref{cl_crazy-combinatorics}(i)
to $A^{(2,3)}$ we get
$$
\begin{aligned}
\Sh^{A^{(2,3)}}_{p,q+r}\circ\Sh^{A''}_{q,r}[p]=\, &
(-1)^{q+r}\cdot\Sh^{D^{(2,3)}}_{p-1,q+r}\circ
A[\one_{[p]},\eta_{q+r},\eta_{q+r}]\circ\Sh^{A''}_{q,r}[p] \\
& +\Sh^{D^{(2,3)}}_{p,q+r-1}\circ A[\eta_p,\one_{[q+r]},\one_{[q+r]}]
\circ\Sh^{A''}_{q,r}[p] \\
=\, & (-1)^{q+r}\cdot\Sh^{D^{(2,3)}}_{p-1,q+r}\circ
\Sh^{D''}_{q,r}[p]\circ A[\one_{[p]},\eta_q,\eta_r] \\
& +\Sh^{D^{(2,3)}}_{p,q+r-1}\circ\Sh^{A''}_{q,r}[p+1]\circ
A[\eta_p,\one_{[q]},\one_{[r]}] &
\quad\text{(by \eqref{eq_with-aligned-for-A})}
\end{aligned}
$$
and after applying again claim \ref{cl_crazy-combinatorics}(i)
to $A''$ and comparing with the foregoing expression for
$\Sh^{A^{(1,2)}}_{p+q,r}\circ\Sh^{A'}_{p,q}[r]$, we obtain the
sought identity.
\end{proof}

\begin{theorem}[Dold-Puppe-Kan]\label{th_Dold-Kan}
For any abelian category $\cA$, and any $k\in\N$,
we have :
\begin{enumerate}
\item
The functors $\sN_\cA$ and $\sN_{\cA,k}$ are equivalences.
\item
If $f,g$ are any two morphisms in $s.\cA$, then there exists
a simplicial homotopy from $f$ to $g$ if and only if there
exists a chain homotopy from $N_\bullet f$ to $N_\bullet g$.
\end{enumerate}
\end{theorem}
\begin{proof} We easily reduce to the case where $\cA$ is
small, and then there exists a fully faithful embedding
$\cA\to\cB$, where $\cB$ is a complete and cocomplete abelian
tensor category, with internal $\Hom$ functor (lemma
\ref{lem_simple-imbeddings}).

(i): We first construct an explicit quasi-inverse for the
functors $\sN_\cB$ and $\sN_{\cB,k}$, as follows. For every
$i\in\N$, consider the cochain complex $\K\La i\Ra_\bullet$
defined in example \ref{ex_simplicial-exponent}(ii); notice
that every morphism $\phi:[i]\to[i']$ in $\Delta$ induces a
morphism
$$
\K\La\phi\Ra_\bullet:=N_\bullet\Z^{\bDelta_\phi}:
\K\La i\Ra_\bullet\to\K\La i'\Ra_\bullet
$$
(notation of remark \ref{rem_cosimpl-simpl}(iii)).
Hence, the system $(\K\La i\Ra_\bullet~|~i\in\N)$ amounts to a cosimplicial
object of $\sC^{\leq 0}(\Z\Mod)$. Now, let $U$ be a unit of the tensor
category $\cB$; for any object $C_\bullet$ of $\sC^{\leq 0}(\cB)$
(resp. of $\sC^{[-k,0]}(\cB)$), we set
$$
K_C[i]:=\cHom_{\sC(\cB)}(\K\La i\Ra_\bullet\otimes_\Z U[0],C_\bullet)
\qquad
\text{for every $i\in\N$ (resp. for every $i\leq k$)}
$$
where $\cHom_{\sC(\cB)}$ is the functor constructed in example
\ref{ex_monoidal}(vii), and the mixed tensor product is defined
as in \eqref{subsec_cplx-mixed-tens}. By the foregoing, it is
clear that the system $(K_C[i]~|~i\in\N)$ amounts to an object
of $s.\cB$ (resp. of $s_k.\cB$). In order to compute $N_\bullet K_C$,
let us set
$$
\Z^{\bDelta_i}_+:=
\sum_{n=1}^i\Img\,(\Z^{\bDelta_{\eps_n}}:\Z^{\bDelta_{i-1}}\to\Z^{\bDelta_i})
\qquad
\Z^{\bDelta_i}_0:=\Z^{\bDelta_i}/\Z^{\bDelta_i}_+
\qquad
\bar\K\La i\Ra_\bullet:=N_\bullet\Z^{\bDelta_i}_0
$$
for every $i>0$, as well as $\Z^{\bDelta_0}_0:=\Z^{\bDelta_0}$
and $\bar\K\La 0\Ra_\bullet:=\K\La 0\Ra_\bullet$; thus,
$\bar\K\La i\Ra_\bullet$ is also the quotient of $\K\La i\Ra_\bullet$
by the sum of the images of the morphisms $\K\La\eps_n\Ra$, for
$n=1,\dots,i$. With this notation, a simple inspection of the
definition shows that
$$
N_iK_C=\cHom_{\sC(\cB)}(\bar\K\La i\Ra_\bullet\otimes_\Z U[0],C_\bullet)
\qquad
\text{for every $i\in\N$ (resp. for every $i\leq k$)}.
$$
On the other hand, using the explicit description of example
\ref{ex_simplicial-exponent}(ii), it is easily seen that
$$
\bar\K\La i\Ra_n\isom\left\{\begin{array}{ll}
                        \Z & \qquad \text{if $n=i$ or $n=i-1\geq 0$} \\
                         0 & \qquad \text{otherwise}.
                       \end{array}\right.
$$
More precisely, $\bar\K\La i\Ra_i$ (resp. $\bar\K\La i\Ra_{i-1}$)
is generated by the basis element $e_\one$ of $\K\La i\Ra_i$ (resp.
$e_{\eps_0}$ of $\K\La i\Ra_{i-1}$) corresponding to the identity
map $\one_{[i]}$ (resp. corresponding to $\eps_0:[i-1]\to[i]$).
Furthermore, example \ref{ex_simplicial-exponent}(i) shows
that the differential $\K\La i\Ra_i\to\K\La i\Ra_{i-1}$ maps $e_\one$
to $e_{\eps_0}$, so it corresponds to the identity map $\Z\to\Z$,
under the foregoing identification. Taking into account
remark \ref{rem_Hom-how-to}(iv), we deduce that $N_iK_C$
is the kernel of the morphism
\set\begin{equation}\label{eq_maps-of-Cs}
C_i\oplus C_{i-1}\to C_{i-1}\oplus C_{i-2}
\end{equation}
given by the matrix
$$
\left(\begin{array}{rr}
(-1)^i\cdot d^C_i & \one_{C_{i-1}} \\
0 & (-1)^{i+1}\cdot d^C_{i-1}
\end{array}\right)
$$
and since $d^C_{i-1}\circ d^C_i=0$, the latter is just $C_i$;
more precisely, $C_i$ is identified with this kernel, via
the monomorphism
\set\begin{equation}\label{eq_identify-C_i-and-Ker}
(\one_{C_i},(-1)^{i+1}\cdot d^C_i):C_i\to C_i\oplus C_{i-1}.
\end{equation}
This identification $N_iK_C\isom C_i$ can also
be described as follows. For every $Z\in\Ob(\cB)$, let
$$
\Hom_\cB(Z,C_\bullet)
$$
be the cochain complex such that
$\Hom_\cB(Z,C_\bullet)_n:=\Hom_\cB(Z,C_n)$ for every $n\in\Z$,
with differentials induced by those of $C_\bullet$, in the
obvious way; then we have natural identifications
$$
\Hom_{\sC(\Z\Mod)}(\bar\K\La i\Ra_\bullet,\Hom_\cB(Z,C_\bullet))\isom
\Hom_{\sC(\cB)}(\bar\K\La i\Ra_\bullet\otimes_\Z Z[0],C_\bullet)\isom
\Hom_\cB(Z,N_iK_C)
$$
(see example \ref{ex_monoidal}(vii)) whose composition with
the induced isomorphism
$$
\Hom_\cB(Z,N_iK_C)\isom\Hom_\cB(Z,C_i)
$$
is given by the rule : 
\set\begin{equation}\label{eq_Dold-explicit}
(\phi_\bullet:\bar\K\La i\Ra_\bullet\to\Hom_\cB(Z,C_\bullet))
\mapsto(\phi_i(\eps_\one):Z\to C_i).
\end{equation}
It remains to determine the differential
$d^N_{i+1}:N_{i+1}K_C\to N_iK_C$; by definition, the latter is
induced by the morphism
$\bar\K\La\eps_0\Ra_\bullet:
\bar\K\La i\Ra_\bullet\to\bar\K\La i+1\Ra_\bullet$.
In turn, the foregoing description says that
$\bar\K\La\eps_0\Ra_\bullet$ is naturally identified
with the morphism of $\sC(\Z\Mod)$
$$
\xymatrix{
\bar\K\La i\Ra_\bullet \ar[d] &
0\ar[r] & 0 \ar[r] \ar[d] & \Z \ar[d]_{\one_\Z} \ar[r] &
\Z \ar[d] \ar[r] & 0 \\
\bar\K\La i+1\Ra_\bullet &
0 \ar[r] & \Z \ar[r] & \Z \ar[r] & 0 \ar[r] & 0
}$$
(where the two copies of $\Z$ on the top horizontal row are
placed in homological degrees $i$ and $i-1$), so $d^N_{i+1}$
is deduced from the morphism in $\sC(\cB)$
$$
\xymatrix{
0 \ar[r] & C_{i+1}\oplus C_i \ar[r] \ar[d] &
C_i\oplus C_{i-1} \ar[d]_\one \ar[r] & 0 \ar[r] \ar[d] & 0 \\
0\ar[r] & 0 \ar[r] & C_i\oplus C_{i-1} \ar[r] &
C_{i-1}\oplus C_{i-2} \ar[r] & 0
}$$
whose rows are given by the morphisms \eqref{eq_maps-of-Cs}.
Therefore, via the identification \eqref{eq_identify-C_i-and-Ker},
the morphism $d^N_{i+1}$ becomes none else than
$(-1)^{i+1}\cdot d^C_{i+1}$, {\em i.e.} we have obtained a natural
isomorphism
$$
\omega^C_\bullet:N_\bullet K_C\isom C_\bullet
$$
for every complex in $\sC^{\leq 0}(\cB)$ (resp. in $\sC^{[-k,0]}(\cB)$).

Conversely, let $B[\bullet]$ be any object of $s.\cB$, and $Z$
any object of $\cB$; clearly
$$
N_\bullet(\bDelta_k\otimes s.Z)=\K\La i\Ra_\bullet\otimes_\Z Z[0]
\qquad
\text{for every $i\in\N$}.
$$
In view of example \ref{ex_monoidal}(vii) and remark
\ref{rem_cosimpl-simpl}(ii), we deduce a natural transformation
\set\begin{equation}\label{eq_Yoneda-by-the-kilo}
{\diagram
\Hom_\cB(Z,B[i]) \ar[r]^-a &
\Hom_{s.\cB}(\bDelta_i\otimes s.Z,B) \ar[d]^b \\
\Hom_\cB(Z,K_{N_\bullet B}[i]) &
\Hom_{\sC(\cB)}(\K\La i\Ra_\bullet\otimes_\Z Z[0],N_\bullet B) \ar[l]_-c
\enddiagram}
\end{equation}
 which, by Yoneda's lemma, comes from a unique morphism in $\cB$
$$
\psi_i^B:B[i]\to K_{N_\bullet B}[i]
\qquad
\text{for every $i\in\N$}.
$$
The same construction applies, in case $B$ is an object of
$s_k.\cB$ : we need only replace the group
$\Hom_{s.\cB}(\bDelta_i\otimes s.Z,B)$ by
$\Hom_{s_k.\cB}(s.\trunc_k(\bDelta_i\otimes s.Z),B)$ in
\eqref{eq_Yoneda-by-the-kilo} : see remark \ref{rem_cosimpl-simpl}(iv).
Moreover, remark \ref{rem_cosimpl-simpl}(iii) implies that
the system $\psi^B:=(\psi^B_i~|~i\in\N)$ amounts to a morphism
$B\to K_{N_\bullet B}$ in $s.\cB$ (resp. in $s_k.\cB$). By the same
token, the $\Z$-linear isomorphism $a$ maps the abelian subgroup
$\Hom_\cB(Z,N_iB)$ isomorphically onto the subgroup
$$
\Hom_{s.\cB}(\Z^{\bDelta_i}_0\otimes_\Z s.Z,B)\isom
\Hom_{s.\Z\Mod}(\Z^{\bDelta_i}_0,\Hom_\cB(Z,B))
$$
for every $Z\in\Ob(\cB)$ and every $i\in\N$, where $\Hom_\cB(Z,B)$
is the simplicial abelian group as in example
\ref{ex_simplicial-exponent}(iii). Explicitly, if
$\beta:Z\to N_iB$ is any morphism, then $a(\beta)$ is the unique
morphism $\Z^{\bDelta_i}_0\to\Hom_\cB(Z,B)$ of simplicial abelian
groups such that $a(\beta)(e_\one)=\beta$, where
$e_\one\in\bar\K\La i\Ra_i\subset\Z^{\bDelta_i}_0[i]$ is the basis
element described in the foregoing. Likewise, we have natural
identifications
$$
\Hom_{\sC(\Z\Mod)}(\K\La i\Ra_\bullet,\Hom_\cB(Z,N_\bullet B))\isom
\Hom_{\sC(\cB)}(\K\La i\Ra_\bullet\otimes_\Z Z[0],N_\bullet B)
$$
and $b$ restricts to a map
$$
\Hom_{s.\Z\Mod}(\Z^{\bDelta_i}_0,\Hom_\cB(Z,B))\to
\Hom_{\sC(\Z\Mod)}(\bar\K\La i\Ra_\bullet,\Hom_\cB(Z,N_\bullet B))
\qquad
\phi\mapsto N_\bullet\phi.
$$
Again, the same applies to an object of $s_k.\cB$, by taking
suitable truncated variants of the above constructions.
Taking into account \eqref{eq_Dold-explicit}, we conclude that
$\psi^B$ induces an isomorphism
$$
N_\bullet\psi^B:N_\bullet B\isom N_\bullet K_{N_\bullet B}
$$
which is inverse to $\omega^{N_\bullet B}$. To finish the proof of
assertion (i) for the category $\cB$, it then suffices to remark :

\begin{claim} The functors $\sN_\cB$ and $\sN_{\cB,k}$ are
conservative.
\end{claim}
\begin{pfclaim} We show, by induction on $n$, that if
a morphism $h$ in $s.\cB$ or $s_k.\cB$ (for any $k\in\N$)
induces an isomorphism $N_\bullet h$, then $h[n]$ is an
isomorphism for every $n\in\N$ (resp. for every $n\leq k$).
The assertion is obvious for $n=0$, hence suppose that $n>0$,
and that $h[n-1]$ is known to be an isomorphism whenever $h$
is a morphism in $s.\cB$ or in $s_k.\cB$ (for arbitrary $k\in\N$)
such that $N_\bullet h$ is an isomorphism. Let $h:A\to B$ be any
such morphism in $s.\cB$ or in $s_k.\cB$. If $h$ is a morphism
in $s_0.\cB$, we are done, hence we may suppose that either $h$
is a morphism in $s.\cB$ or $k>0$. Set
$$
A':=\Ker(g_A:\gamma A\to A)
\qquad
\text{(\ resp. $A':=\Ker(g_A:\gamma_{k-1}A\to s.\trunc_{k-1}A)$\ )}
$$
as well as $B':=\Ker(g_B)$ (notation of remark
\ref{rem_path-spaces}(ii)). Notice that $g_A$ is an epimorphism,
since $\partial_{n+2}$ admits the section $\sigma_{n+1}$, for
every $n\in\N$ (resp. for every $n\leq k$). Therefore, we have
a commutative diagram in $s.\cB$ with exact rows 
$$
\xymatrix{
0\ar[r] & A' \ar[r] \ar[d]_{h'} &
\gamma A \ar[r] \ar[d]_{\gamma h} & A \ar[r] \ar[d]^h & 0 \\
0\ar[r] & A' \ar[r] & \gamma A \ar[r] & A \ar[r] & 0
}$$
(resp. a corresponding diagram in $s_{k-1}.\cB$).
By inspecting the definitions, it is easily seen that
$N_\bullet A'=(N_\bullet A)[-1]$, $N_\bullet A'=(N_\bullet A)[-1]$,
and $N_\bullet h'=(N_\bullet h)[-1]$; especially, $N_\bullet h'$
is an isomorphism, so $h'[n-1]$ is an isomorphism, by inductive
assumption. The same holds also for $h[n-1]$, and we conclude
that $\gamma h[n]$ is an isomorphism. But $\gamma h[n]=h[n+1]$,
so we are done.
\end{pfclaim}

Lastly, in order to prove assertion (i) for the original category
$\cA$, it suffices to notice :

\begin{claim} Let $C_\bullet$ be any object of $\sC(\cA)$ (resp.
of $\sC^{[-k,0]}(\cA)$), and regard $C_\bullet$ as an object of
$\sC(\cB)$ (resp. of $\sC^{[-k,0]}(\cB)$), via the fully faithful
embedding $\cA\to\cB$. Then $K_C$ is isomorphic to an object of
$s.\cA$ (resp. $s_k.\cA$), regarded as a full subcategory of
$s.\cB$ (resp. of $s_k.\cB$), via the same embedding.
\end{claim}
\begin{pfclaim} This follows easily, by remarking that
$\K\La i\Ra_\bullet$ lies in $\sC^{[-i,0]}(\Z\Mod)$ for every $i\in\N$,
and $\K\La i\Ra_j$ is a finitely generated abelian group for every
$i,j\in\N$ : details left to the reader.
\end{pfclaim}

(ii): First, let $f,g:A\to B$ be two morphisms in $s.\cA$,
and $u:\bDelta_1\otimes A\to B$ a homotopy from $f$ to $g$
(see remark \ref{rem_mixed-simpl-tensors}(v)); especially,
$$
u\circ(\bDelta_{\eps_1}\otimes A)=f
\qquad\text{and}\qquad
u\circ(\bDelta_{\eps_0}\otimes A)=g.
$$
Notice that
$$
\Z^{\bDelta_1}_\bullet\otimes_\Z A_\bullet=
\Tot(\bDelta_1\boxtimes A)_{\bullet\bullet}
$$
(notation of remark \ref{rem_mixed-simpl-tensors}(iii) and
example \ref{ex_simplicial-exponent}(i)); there follows a
morphism in $\sC(\cA)$
$$
\tilde u_\bullet:
\Z^{\bDelta_1}_\bullet\otimes_\Z A_\bullet\xrightarrow{\ \Sh_\bullet\ }
(\bDelta_1\otimes A)_\bullet\xrightarrow{\ \sq_\bullet\ }
N_\bullet(\bDelta_1\otimes A)\xrightarrow{\ N_\bullet u\ }N_\bullet B
$$
where $\Sh_\bullet$ denotes the shuffle map for the bisimplicial
object $A\boxtimes\bDelta_1$, and $\sq_\bullet$ is the projection
defined in \eqref{subsec_sj-and-sq}. Moreover, the maps
$\bDelta_{\eps_i}[0]:\bDelta_0[0]\to\bDelta_1[0]$ ($i=0,1$)
induce morphisms
$$
\tilde e_{i,n}:=\Z^{\bDelta_{\eps_i}[0]}\otimes_\Z A_n:
A_n\to\Z^{\bDelta_1[0]}\otimes_\Z A_n\subset
(\Z^{\bDelta_1}_\bullet\otimes_\Z A_\bullet)_n
$$
that amount to morphisms of cochain complexes
$\tilde e_{i,\bullet}:A_\bullet\to\Z^{\bDelta_1}_\bullet\otimes_\Z A_\bullet$
($i=0,1$), and a simple inspection of the definitions shows
that
$$
\Sh_\bullet\circ\tilde e_{i,\bullet}=(\Delta_{\eps_i}\otimes A)_\bullet
\qquad
\text{for $i=0,1$}
$$
whence
$$
\tilde u_\bullet\circ\tilde e_{1,\bullet}=f_\bullet
\qquad\text{and}\qquad
\tilde u_\bullet\circ\tilde e_{0,\bullet}=g_\bullet.
$$
The construction makes it clear that $\tilde e_i$ restricts
to a morphism $N_\bullet A\to\K\La 1\Ra_\bullet\otimes_\Z N_\bullet A$,
and the latter is none else than the map
$\iota_i\otimes_\Z N_\bullet A$, with the notation of remark
\ref{rem_chain-homotopies}(ii). We conclude that the morphism
$$
\bar u_\bullet:\K\La 1\Ra_\bullet\otimes_\Z N_\bullet A\hookrightarrow
\Z^{\bDelta_1}_\bullet\otimes_\Z A_\bullet\xrightarrow{\ \tilde u_\bullet\ }
N_\bullet B
$$
is a homotopy from $N_\bullet f$ to $N_\bullet g$ (notation of
\eqref{subsec_sj-and-sq}).

Conversely, suppose that
$\bar v_\bullet:\K\La 1\Ra_\bullet\otimes_\Z N_\bullet A\to N_\bullet B$
is a homotopy from $N_\bullet f$ to $N_\bullet g$. Since
$\K\La 1\Ra_\bullet$ is a direct summand of $\Z^{\bDelta_1}_\bullet$
and $N_\bullet A$ is a direct summand of $A_\bullet$, we may
extend $\bar v_\bullet$ to a morphism
$$
\tilde v_\bullet:\Z^{\bDelta_1}_\bullet\otimes_\Z A_\bullet\to N_\bullet B
$$
such that $\tilde v_\bullet$ is the zero morphism on the direct summands
other than $\K\La 1\Ra_\bullet\otimes_\Z N_\bullet A$. Next, consider
the composition
$$
v_\bullet:N_\bullet(\bDelta_1\otimes A)\xrightarrow{\ \sj_\bullet\ }
(\bDelta_1\otimes A)_\bullet\xrightarrow{\ \AW_\bullet\ }
A_\bullet\otimes_\Z\Z^{\bDelta_1}_\bullet\xrightarrow{\ \Psi_\bullet\ }
\Z^{\bDelta_1}_\bullet\otimes_\Z A_\bullet\xrightarrow{\ \tilde v_\bullet\ }
N_\bullet B
$$
where $\sj_\bullet$ is the natural injection (see
\eqref{subsec_sj-and-sq}), $\AW_\bullet$ is the Alexander-Whitney
map for the bisimplicial object $A\boxtimes\bDelta_1$ (notice that
$\bDelta_1\otimes A=(A\boxtimes\bDelta_1)^\Delta$), and
$\Psi_\bullet$ is the commutativity constraint (see example
\ref{ex_monoidal}(i)). By (i), the morphism $v_\bullet$ comes
from a unique morphism
$$
v:\bDelta_1\otimes A\to B
\qquad
\text{in $s.\cA$}.
$$
On the other hand, since $N_{p+q}A$ is contained in the kernel
of $A[\eps^{q\vee}_{p,0}]$ for every $p,q\in\N$, it is easily seen
that the diagram
$$
\xymatrix{
N_\bullet A \ar[rr]^-{\iota_i\otimes_\Z N_\bullet A}
\ar[d]_{N_\bullet(\bDelta_{\eps_i}\otimes A)}
& & \K\La 1\Ra_\bullet\otimes_\Z N_\bullet A \ar[d] \\
N_\bullet(\bDelta_1\otimes A) 
\ar[rr]^-{\Psi_\bullet\circ\AW_\bullet\circ\sj_\bullet} & &
\Z^{\bDelta_1}_\bullet\otimes_\Z A_\bullet
}$$
commutes for $i=0,1$. We conclude that $v$ is a homotopy
from $f$ to $g$.
\end{proof}

\begin{corollary}\label{cor_Dold-Kan}
Let $\cA$ be any abelian category. We have :
\begin{enumerate}
\item
If $k\in\N$ is any integer, and $A$ any $k$-truncated simplicial
object of $\cA$, then
$$
H_i(\cosk_kA)=0
\qquad
\text{for every $i\geq k$}. 
$$
\item
Every homotopically trivial augmented simplicial object of
$\cA$ is aspherical.
\end{enumerate}
\end{corollary}
\begin{proof}(i): Denote by
$t_{\leq k}:\sC^{\leq 0}(\cA)\to\sC^{[-k,0]}(\cA)$ the brutal
truncation functor (see \eqref{subsec_brutal-truncate}). In
light of theorem \ref{th_Dold-Kan}, we see that $t_{\leq k}$
admits a right adjoint $v_k:\sC^{[-k,0]}(\cA)\to\sC^{\leq 0}(\cA)$,
and clearly there are natural isomorphisms
$$
N_\bullet\cosk_kA\isom v_k N_\bullet A
\qquad
\text{for every $A\in\Ob(s_k.\cA)$}.
$$
Taking into account theorem \ref{th_degen-plus-norm}(iii),
we are then reduced to showing

\begin{claim} $H_i(v_kK_\bullet)=0$ for every
$(K_\bullet,d_\bullet)\in\Ob(\sC^{[-k,0]})$ and every $i\geq k$.
\end{claim}
\begin{pfclaim} Indeed it is easily seen that :
$$
(v_kK_\bullet)_i=\left\{
\begin{array}{ll}K_i & \text{for $\leq k$} \\
           \Ker\,d_k & \text{for $i=k+1$} \\
                   0 & \text{for $i>k+1$}
\end{array}\right.
$$
and the differential of $v_kK_\bullet$ in degree $\leq k$
agrees with that of $K_\bullet$, whereas in degree $k+1$
it is the natural inclusion map (details left to the reader).
The claim follows immediately.
\end{pfclaim}

(ii) follows directly from theorems \ref{th_Dold-Kan}(ii)
and \ref{th_degen-plus-norm}(iii), and remark
\ref{rem_abel-homotopies}(ii).
\end{proof}

\subsection{Simplicial sets}\label{sec_simplicial-sets}
This section, which complements the previous one, collects
some classical material pertaining to the homotopy theory of
the category of simplicial sets. The presentation is borrowed
from \cite{Go-Ja} and \cite{Joy}, where much more may be found.

\sset\subsubsection{}\label{subsec_eff-sset-epis}
To begin with, notice that the category $s.\Set$ of
{\em simplicial sets} is none else than the category
of presheaves on the category $\Delta$, so $s.\Set$
is complete and cocomplete, and all limits and colimits
in $s.\Set$ are computed argumentwise (see
\eqref{sec_topoi}); also, all colimits and monomorphisms
are universal, and all epimorphisms are universal effective.

\begin{example}\label{ex_simplicial-sets}
(i)\ \
We have already introduced in example \ref{ex_simplicies}
the simplicial set $\bDelta_k$ (for any $k\in\N$), which
is just the image of $[k]$ under the Yoneda embedding
$$
h:\Delta\to s.\Set.
$$
To ease notation, for any morphism $\phi:[n]\to[m]$ in
$\Delta$, we shall usually write
$\phi:\bDelta_n\to\bDelta_m$ instead of $\bDelta_\phi$.
By Yoneda's lemma (proposition \ref{prop_yoneda}), we
have natural identifications
\set\begin{equation}\label{eq_yon-s.sets}
A[n]\isom\Hom_{s.\Set}(\bDelta_n,A)
\qquad
\text{for every $A\in\Ob(s.\Set)$ and every $n\in\N$}.
\end{equation}

(ii)\ \ 
More generally, the category of simplicial sets contains
all the products $\bDelta_n\times\bDelta_m$ (for every
$m,n\in\N$). Furthermore, if $S$ is any set and $A$ any
simplicial set, we may form the product $s.S\times A$
(notation of \eqref{subsec_simplicial-object}), which
we shall denote simply by $S\times A$. Explicitly, we have
$$
(S\times A)[n]:=S\times A[n]
\qquad
\text{for every $n\in\N$}
$$
and the faces and degeneracies of $S\times A$ are derived
from those of $A$, in the obvious fashion.

(iii)\ \
Moreover, lemma \ref{lem_lable} yields a natural
presentation
$$
A\isom\colim_{h\Delta/A}h_{s.\Set}\circ\iota_A
$$
for every simplicial set $A$. By inspecting the definitions,
we see that the latter amounts to the following description
of $A$. Notice the natural identification
$\bDelta_j[i]\isom\Hom_{s.\Set}(\bDelta_i,\bDelta_j)$ provided
by \eqref{eq_yon-s.sets}; we have a natural diagram in $s.\Set$
\set\begin{equation}\label{eq_attach-cells}
\coprod_{i,j\in\N}(A[i]\times\bDelta_i[j])\times\bDelta_j
\xymatrix{\ar@<.5ex>[r]^-{p_\bullet} \ar@<-.5ex>[r]_-{c_\bullet} & }
\coprod_{n\in\N} A[n]\times\bDelta_n \xrightarrow{\ t_\bullet\ }A
\end{equation}
where $p_\bullet$ are $c_\bullet$ are the (unique) morphisms
which restrict respectively to morphisms
$$
\{a\}\times\bDelta_j\xleftarrow{\ \phi\ }
\{(a,\phi)\}\times\bDelta_i\xrightarrow{\ \one_{\bDelta_i}\ }
\{A[\phi](a)\}\times\bDelta_i
$$
for every $i,j\in\N$ and every
$(a,\phi)\in A[i]\times\bDelta_j[i]$, and where $t_\bullet$
is the (unique) morphism whose restriction to 
$\{b\}\times\bDelta_n\to A$ is the morphism corresponding
to $b\in A[n]$ under the identification \eqref{eq_yon-s.sets},
for every $b\in A[n]$.
Then $t_\bullet$ is an epimorphism, and \eqref{eq_attach-cells}
naturally identifies $A$ with the coequalizer of $p_\bullet$
and $c_\bullet$.

(iv)\ \
Another useful simplicial set is the {\em boundary} of
$\bDelta_k$, denoted
$$
\partial\bDelta_k
\qquad
\text{for every integer $k>0$}
$$
which is defined as the smallest subobject of $\bDelta_k$
containing the images of all the face morphisms
$\eps_i:\bDelta_{k-1}\to\bDelta_k$, for $i=0,\dots,k$
(see \eqref{subsec_do-faces}). Thus, we have a natural
epimorphism in $s.\Set$
\set\begin{equation}\label{eq_present-boundary}
[k]\times\bDelta_{k-1}\to\partial\bDelta_k
\qquad
\text{for every $k>0$}
\end{equation}
whose restriction to the subobject $\{i\}\times\bDelta_{k-1}$
agrees with $\eps_i$, for every $i=0,\dots,k$.

Notice that the simplicial identities yield a cartesian
diagram in $\Delta^\wedge$
$$
{\diagram [k-2] \ar[r]^-{\eps_{j-1}} \ar[d]_{\eps_i} &
          [k-1] \ar[d]^{\eps_i} \\
          [k-1] \ar[r]^-{\eps_j} & [k]
\enddiagram}
\qquad
\text{whenever $0\leq i<j\leq k$}.
$$
Since the Yoneda embedding commutes with representable limits
(corollary \ref{cor_pre-misc}(vi)), there follows a cartesian
diagram in $s.\Set$
$$
{\diagram \bDelta_{k-2} \ar[r]^-{\eps_{j-1}} \ar[d]_{\eps_i} &
          \bDelta_{k-1} \ar[d]^{\eps_i} \\
          \bDelta_{k-1} \ar[r]^-{\eps_j} & \bDelta_k
\enddiagram}
\qquad
\text{for $0\leq i<j\leq k$}
$$
(details left to the reader). Now, denote by
$S_k\subset[k]\times[k]$ the subset of all pairs $(i,j)$
with $i<j$; since \eqref{eq_present-boundary} is effective
(by \eqref{subsec_eff-sset-epis}), we deduce a commutative
diagram
$$
\xymatrix{S_k\times\bDelta_{k-2}
\ar@<.5ex>[r]^-{\eps'_\bullet} \ar@<-.5ex>[r]_-{\eps''_\bullet}
& [k]\times\bDelta_{k-1} \ar[r] & \partial\bDelta_k
}\qquad
\text{in $s.\Set$}
$$
which presents $\partial\bDelta_k$ as the coequalizer of
$\eps'_\bullet$ and $\eps''_\bullet$, where $\eps'_\bullet$
(resp. $\eps''_\bullet$) is the morphism whose restriction
to $\{(i,j)\}\times\bDelta_{k-2}$ agrees with
$\eps_i:\{(i,j)\}\times\bDelta_{k-2}\to\{j\}\times\bDelta_{k-1}$
(resp. with
$\eps_{j-1}:\{(i,j)\}\times\bDelta_{k-2}\to\{i\}\times\bDelta_{k-1}$)
for every $(i,j)\in S_k$. As an immediate corollary, we see that
if $A$ is any simplicial set, a morphism $\partial\bDelta_k\to A$
is the same as the datum of an ordered sequence
$x_0,\dots,x_k$ of $k+1$ elements of $A[k-1]$, such that
$$
\partial_ix_j=\partial_{j-1}x_i
\qquad
\text{whenever $0\leq i<j\leq k$}.
$$
Moreover, if $x\in A[n]$ is any element, we may regard $x$
as a morphism $\bar x:\bDelta_n\to A$ via the natural
identification \eqref{eq_yon-s.sets}, and then we shall
often denote by
$$
\partial x:\partial\bDelta_n\to A
$$
the restriction of $\bar x$ to the suboboject
$\partial\bDelta_n$. Lastly, we set
$\partial\bDelta_0:=s.\emptyset$.

(v)\ \
Another important simplicial set is $\bLambda_k^n$, which is
defined for every $k,n\in\N$ such that $k>0$ and $n=0,\dots,k$;
namely, it is the smallest subobject of $\bDelta_k$ that
contains the images of all the face morphisms
$\eps_i:\Delta_{k-1}\to\Delta_k$, except for the face $\eps_n$.
Thus, for every such $k$ and $n$ we have a natural monomorphism
$$
\iota_k^n:\bLambda_k^n\to\bDelta_k.
$$
The same argument as in (iv) yields a presentation of
$\bLambda^n_k$ as the coequalizer of two morphisms :
$$
\xymatrix{S^n_k\times\bDelta_{k-2}
\ar@<.5ex>[r]^-{\eps'_\bullet} \ar@<-.5ex>[r]_-{\eps''_\bullet}
& ([k]\setminus\{n\})\times\bDelta_{k-1} \ar[r] & \bLambda^n_k
}$$
where $S^n_k\subset S_k$ is the subset of all pairs
$(i,j)$ with $i,j\neq n$, and $\eps'_\bullet$, $\eps''_\bullet$
are the restrictions of the morphisms with the same name
appearing in (iii).

(vi)\ \
Let $\psi:X\to Y$ be any morphism of simplicial sets,
$y\in Y[0]$ any element, and $j_y:\Delta_0\to Y$ the
corresponding morphism of simplicial sets; the resulting
fibre product
$$
\psi^{-1}(y):=\Delta_0\times_YX
$$
is called the {\em fibre of\/ $\psi$ over $y$}. Explicitly,
set $\{y_n\}:=\Img\,j_y[n]$ for every $n\in\N$. Then
$$
\psi^{-1}(y)[n]=\psi[n]^{-1}(y_n)
\qquad
\text{for every $n\in\N$}
$$
and the faces and degeneracies of $\psi^{-1}(y)$ are
the restrictions of the corresponding maps for $X$.
\end{example}

\begin{remark}\label{rem_degenerate}
(i)\ \
From the simplicial identities of \eqref{subsec_do-faces}
we also get a commutative diagram
$$
\xymatrix{
\bDelta_{n+1} \ar[r]^-{\eta_i} \ar[d]_{\eta_{j+1}} &
\bDelta_n \ar[d]^{\eta_j} \\
\bDelta_n \ar[r]^-{\eta_i} & \bDelta_{n-1}
}$$
for every integer $n>0$ and every $i,j\leq n-1$ with $i\leq j$.
We claim that this diagram is cocartesian in $s.\Set$.
Indeed, let $X$ be any simplicial set, and $f,g:\bDelta_n\to X$
any two morphisms such that $f\circ\eta_i=g\circ\eta_{j+1}$.
We set $h:=f\circ\eps_{j+1}:\bDelta_{n-1}\to X$, and we notice
that
$$
h=f\circ\eta_i\circ\eps_i\circ\eps_{j+1}=
g\circ\eta_{j+1}\circ\eps_i\circ\eps_{j+1}=
g\circ\eta_{j+1}\circ\eps_{j+2}\circ\eps_i=
g\circ\eps_i.
$$
Therefore :
$$
\begin{aligned}
h\circ\eta_j=\, & g\circ\eps_i\circ\eta_j=
g\circ\eta_{j+1}\circ\eps_i=f\circ\eta_i\circ\eps_i=f \\
h\circ\eta_i=\, & f\circ\eps_{j+1}\circ\eta_i=
f\circ\eta_i\circ\eps_{j+2}=
g\circ\eta_{j+1}\circ\eps_{j+2}=g.
\end{aligned}
$$
Lastly, since both $\eta_i$ and $\eta_j$ are epimorphisms,
$h$ is uniquely determined by either of the identities
$h\circ\eta_j=f$ and $h\circ\eta_i=g$, whence the
assertion.

(ii)\ \
We point out that also the diagram
$$
\xymatrix{
\bDelta_{n+1} \ar[r]^-{\eta_i} \ar[d]_{\eta_i} &
\bDelta_n \ar[d]^{\one_{\bDelta_n}} \\
\bDelta_n \ar[r]^-{\one_{\bDelta_n}} & \bDelta_n
}$$
is trivially cocartesian for every $n\in\N$ and
every $i\leq n$, since $\eta_i$ is an epimorphism.

(iii)\ \
Let $k\in\N$ be any integer. The categories $\Set$ of all
small sets and $\fSet$ of finite sets, both satisfy the
conditions of \eqref{subsec_coskeleton}, so we get left
and right adjoints
$$
\sk_k,\cosk_k:s_k.\Set\to s.\Set
\qquad
\sk_k,\cosk_k:s_k.\fSet\to s.\fSet
$$
for the respective $k$-truncation functors.

(iv)\ \
For any simplicial set $X$ and any integer $r>0$, let
us say that an element $x\in X[r]$ is a
{\em degenerate simplex}, if $x=\sigma_iy$ for some
$i=0,\dots,r-1$ and some $y\in X[r-1]$. We claim that,
for any $k$-truncated simplicial set $Y$, and every
$r>k$, every element of $\sk_kY[r]$ is a degenerate
simplex. Indeed, let $\alpha:[k]\to[r]$ be any morphism
in $\Delta^o$; notice first that, under the identification
$Y[k]\isom\sk_kY[k]$ given by the unit of adjunction, the
map $\sk_kY[\alpha]:\sk_kY[k]\to\sk_kY[r]$ corresponds to
the natural map $j_\alpha:F[k]\to\sk_kY[r]$. But $\sk_kY[r]$
is the union of the images of all such maps (see example
\ref{ex_complete-cats}(i)), and on the other hand, any
such map factors through some degeneracy map $\sigma_i$,
whence the assertion.

(v)\ \
For any simplicial set $X$ and every $k,n\in\N$ with
$n\geq k$, the counits of adjunction give a natural
commutative diagram
$$
\xymatrix{
\sk_k(s.\trunc_kX) \ar[rr]^-{\eps^{(k,n)}_X}
\ar[rd]_{\eps^{(k)}_X} & &
\sk_n(s.\trunc_nX) \ar[ld]^{\eps^{(n)}_X} \\
& X
}$$
amounting to a cocone with vertex $X$, and it follows
easily from the discussion of \eqref{subsec_coskeleton}
that $s.\trunc_k(\eps^{(k,n)}_X)$ is an isomorphism for every
such $k,n\in\N$, and the resulting morphism
$$
\colim_{n\in\N}\sk_n(s.\trunc_nX)\to X
$$
is an isomorphism. We define the {\em $n$-th skeleton of $X$}
as the simplicial subset
$$
\Sk_nX:=\Img\,\eps_X^{(n)}
\qquad
\text{for every $n\in\N$}.
$$
\end{remark}

\begin{lemma}\label{lem_degenerate}
Let $X$ be any simplicial set, $r>0$ any integer, and
$x,y\in X[r]$ be any two degenerate simplices (see
remark {\em\ref{rem_degenerate}(iii)}). The following
holds :
\begin{enumerate}
\item
If $\partial x=\partial y$, then $x=y$.
\item
There exists a unique pair $(p,z)$ where $p:[r]\to[n]$ is
an epimorphism and $z\in X[n]$ is a non-degenerate simplex
such that $x=X[p](z)$.
\item
The natural morphism $\sk_n(s.\trunc_nX)\to\Sk_nX$ is an
isomorphism for every $n\in\N$.
\end{enumerate}
\end{lemma}
\begin{proof}(i): Say that $x=\sigma_mz$ and $y=\sigma_nw$
for some $z,w\in X[r-1]$. If $m=n$, we have
$$
z=\partial_m\sigma_mz=\partial_mx=\partial_my=
\partial_m\sigma_mw=w
$$
whence the assertion. Hence, we may assume that $m<n$,
in which case we get :
$$
z=\partial_mx=\partial_m\sigma_nw=\sigma_{n-1}\partial_mw
$$
so that
$x=\sigma_m\sigma_{n-1}\partial_mw=\sigma_n\sigma_m\partial_mw$.
Thus, we may replace $z$ by $\sigma_m\partial_mw$ and $m$
by $n$, and reduce to the foregoing case, so the proof
is complete.

(ii): The existence of such a pair $(p,z)$ is obvious.
Next, suppose that $(p':[r]\to[n'],z')$ is another such
pair; from remark \ref{rem_degenerate}, it follows easily
that there exists a cocartesian diagram of the form
$$
\xymatrix{
\bDelta_r \ar[r]^-{\bDelta_p} \ar[d]_{\bDelta_{p'}} &
\bDelta_n \ar[d]^{\bDelta_{q'}} \\
\bDelta_{n'} \ar[r]^-{\bDelta_q} & \bDelta_k
}$$
for a suitable $k\leq n,n'$ and epimorphisms $q:[n]\to[k]$
and $q':[n']\to[k]$ (details left to the reader).
Regarding $z$ and $z'$ as morphisms $x:\bDelta_n\to X$
and respectively $x':[n']\to X$, our assumption gives
the identity :
$$
z\circ\bDelta_p=z'\circ\bDelta_{p'}
$$
whence a unique $w\in X[k]$ such that $X[q'](w)=z$ and
$X[q](w)=z'$. Since $z$ and $z'$ are non-degenerate, it
follows that $n=k=n'$ and $q=q'=\one_{\bDelta_k}$, whence
$z=z'$, as stated.

(iii): Set $Y:=\sk_n(s.\trunc_nX)$ for every $n\in\N$;
the natural map $f[k]:Y[k]\to\Sk_n[k]$ is surjective
for every $k\in\N$, and is bijective for every $k\leq n$,
so it suffices to check that $f[k]$ is injective for every
$k>n$. Thus, fix $k>N$, and let $y,y'\in Y[k]$ be any
two simplices such that $f[k](y)=f[k](y')$. Remark
\ref{rem_degenerate}(iv) tells us that $y$ and $y'$ are
degenerate, say $y=Y[p](z)$ and $y'=Y[p'](z')$ for
epimorphisms $p:[k]\to[m]$ and $p:[k]\to[m']$, and
non-degenerate simplices $z\in Y[m]$, $z'\in Y[m']$.
We must then have $m,m'\leq n$, and therefore
$x:=f[m](z)$ and $x':=f[m'](z')$ are two non-degenerate
simplices of $\Sk_nX$, such that
$$
\Sk_n[p](x)=f[k](y)=f[k](y')=\Sk_n[p'](x').
$$
By (ii), we deduce that $p=p'$ and $x=x'$, whence
$y=y'$, as required.
\end{proof}

\begin{remark}\label{rem_skeleta}
Let $X$ be any simplicial set and $n\in\N$ any integer.

(i)\ \
Clearly, $\Sk_nX\subset\Sk_{n+1}X$, and we have
\set\begin{equation}\label{eq_isom-trunc}
X=\bigcup_{r\in\N}\Sk_rX.
\end{equation}
We say that $X$ has dimension $\leq n$, if $X=\Sk_nX$.
Let $s.\Set_{\leq n}$ be the full subcategory of $s.\Set$
whose objects are the simplicial set of dimension $\leq n$;
it follows easily from lemma \ref{lem_degenerate}(iii)
that the adjunction $(\sk_n,\trunc_n)$ establishes an
equivalence
$$
s.\Set_{\leq n}\isom s_k.\Set.
$$

(ii)\ \
Moreover, $\Sk_nX$ can be obtained by 
``attaching $n$-cells'' to $\Sk_{n-1}X$, if $n>0$.
Namely, let $E^n_X\subset X[n]$ denote the set of
non-degenerate $n$-simplices of $X$; we get a
commutative diagram
$$
\cD
\quad :\quad
{\diagram
E^n_X\times\partial\bDelta_n
\ar[rr]^-{E^n_X\times i_n} \ar[d] & &
E^n_X\times\bDelta_n \ar[d]^c \\
\Sk_{n-1}X \ar[rr] & & \Sk_nX
\enddiagram}
$$
where $i_n:\partial\bDelta_n\to\bDelta_n$ and the bottom
horizontal arrow are the natural inclusion maps, and
$c$ is the unique morphism whose restriction to the
factor $\{x\}\times\bDelta_n$ is the morphism
$\bar x:\bDelta_n\to X$ corresponding to $x$, for every
$x\in E^n_X$. We claim that $\cD$ is a cocartesian
diagram. For the proof, notice that all the simplicial
sets in $\cD$ have dimension $\leq n$; by (i), it then
suffices to check that $\trunc_n\cD$ is cocartesian, {\em i.e.}
that the diagram of sets $\cD[k]$ is cocartesian, for
every $k\leq n$ (see \eqref{subsec_eff-sset-epis}).
The latter assertion is clear for $k<n$, since in this
case both horizontal arrows of $\cD[k]$ are isomorphisms.
For $k=n$, notice that the complement $S$ of
$E^n_X\times\partial\bDelta_n[n]$ in $E^n_X\times\bDelta_n[n]$
is naturally identified with $E^n_X$, which is also the
complement of $\Sk_{n-1}X[n]$ in $\Sk_nX[n]$; it then suffices
to remark that, under these natural bijections, the
restriction of $c[n]$ to $S$ is identified with the identity
map $E^n_X\isom E^n_X$ (details left to the reader).
\end{remark}

\begin{definition}\label{def_fibrant-Kan}
Let $\phi:A\to B$ and $\psi:X\to Y$ be two morphisms of
simplicial sets.
\begin{enumerate}
\item
We say that {\em $\psi$ has the right lifting property
with respect to $\phi$}, if for every commutative diagram
$$
{\diagram A \ar[r]^-\alpha \ar[d]_\phi & X \ar[d]^\psi \\
          B \ar[r]^-\beta & Y
\enddiagram}
\qquad
\text{in $s.\Set$}
$$
there exists a morphism $\gamma:B\to X$ such that
$\gamma\circ\phi=\alpha$ and $\psi\circ\gamma=\beta$.
\item
We say that $\psi:X\to Y$ is a {\em Kan fibration} (or
briefly, that $\psi$ is a {\em fibration}), if $\psi$
has the right lifting property with respect to all the
monomorphisms $\iota_k^n:\bLambda^n_k\to\bDelta_k$, for
every $k,n\in\N$ with $k>0$ and $n\leq k$ (notation of
example \ref{ex_simplicial-sets}(v)).
\item
We say that $\phi$ is a {\em retract of $\psi$} if there
exists a commutative diagram of simplicial sets
$$
\xymatrix{ A \ar[r] \ar[d]_\phi \ar@/^1pc/[rr]^-{\one_A} &
           X \ar[d]^\psi \ar[r] &
           A \ar[d]^\phi \\
           B \ar[r] \ar@/_1pc/[rr]_-{\one_B} & Y \ar[r] & B.
}$$
\item
We say that a simplicial set $A$ is {\em fibrant} if the
(unique) morphism $A\to\bDelta_0$ is a fibration.
\end{enumerate}
\end{definition}

\begin{remark}\label{rem_fibrant-Kan}
In light of example \ref{ex_simplicial-sets}(iv,v), it
is clear that a simplicial set $A$ is fibrant if and only
if the following {\em Kan extension condition} holds. For
every $n,k\in\N$ with $0\leq k\leq n+1$, and every sequence
$x_0,\dots,x_{k-1},x_{k+1},\dots,x_{n+1}$ of elements of
$A[n]$ such that
$$
\partial_ix_j=\partial_{j-1}x_i
\qquad
\text{whenever $0\leq i<j\leq n+1$ and $i,j\neq k$}
$$
there exists an element $x\in A[n+1]$ such that
$$
\partial_ix=x_i
\qquad\text{for every $i=0,\dots,k-1,k+1,\dots,n+1$}.
$$
\end{remark}

\begin{example} Let $A$ be any {\em simplicial group},
{\em i.e.} an object of $s.\Grp$, where $\Grp$
denotes the category of groups. Then the simplicial set
underlying $A$ is fibrant. Indeed, suppose that $n$, $k$
and $x_0,\dots,x_{n+1}$ are as in remark
\ref{rem_fibrant-Kan}. We construct by induction on
$i=-1,\dots,n+1$ an element $y_i\in A[n+1]$ such that
$\partial_jy_i=x_j$ whenever $0\leq j\leq i$, $j\neq k$.
Indeed, the condition is fulfilled trivially for $i=-1$,
by setting $y_{-1}:=1$, the neutral element of the group
$A[n+1]$. Suppose that $i\geq 0$, and $y_{i-1}$ has already
been exhibited as required; if $i=k$, we set $y_i:=y_{i-1}$,
which again trivially fulfills the stated condition.
Otherwise, consider the element
$u:=x_i^{-1}\cdot(\partial_iy_{i-1})$, and notice that
$$
\partial_ju=
(\partial_jx_i^{-1})\cdot(\partial_j\partial_iy_{i-1})=
(\partial_jx_i^{-1})\cdot(\partial_{i-1}\partial_jy_{i-1})=
(\partial_jx_i^{-1})\cdot(\partial_{i-1}x_j)=1
$$
for every $j=0,\dots,i-1$ with $j\neq k$. Hence, set
$y_i:=y_{i-1}\cdot(\sigma_iu)^{-1}$; we have
$$
\partial_jy_i=(\partial_jy_{i-1})\cdot(\partial_j\sigma_iu)^{-1}
=x_i\cdot(\sigma_{i-1}\partial_ju)=x_i
\qquad
\text{whenever $0\leq j<i$}.
$$
Likewise :
$$
\partial_iy_i=(\partial_iy_{i-1})\cdot(\partial_i\sigma_iu)^{-1}
=(\partial_iy_{i-1})\cdot(\partial_i\sigma_i\partial_iy_{i-1})^{-1}
\cdot(\partial_i\sigma_ix_i)=x_i
$$
as needed. Thus, the element $x:=y_{n+1}$ fulfills the
condition of remark \ref{rem_fibrant-Kan}.
\end{example}

For any given morphism $\psi$ of $s.\Set$, the set of all
morphisms $\phi$ such that $\psi$ has the right lifting
property with respect to $\phi$ is closed under certain
elementary operations, that are singled out in the following :

\begin{definition}
Let $\Sigma$ be a set of monomorphisms of $s.\Set$.

(i)\ \
We say that $\Sigma$ is {\em saturated} if the following
conditions hold :
\begin{enumerate}
\alphaenu
\item
All isomorphisms of $s.\Set$ lie in $\Sigma$.
\item
For every countable system of morphisms in $s.\Set$
$$
A_0\xrightarrow{\ \phi_0\ }A_1\xrightarrow{\ \phi_1\ }A_2
\to\cdots
$$
such that $\phi_n\in\Sigma$ for every $n\in\N$, also
the induced morphism
$$
A_0\to\colim_{n\in\N}A_n
$$
lies in $\Sigma$.
\item
For any cocartesian diagram of simplicial sets :
$$
\xymatrix{ A \ar[r] \ar[d]_\phi & A' \ar[d]^{\phi'} \\
           B \ar[r] & B\amalg_AA'
}$$
such that $\phi\in\Sigma$, we have $\phi'\in\Sigma$.
\item
If $(\phi_i:A_i\to B_i~|~i\in I)$ is an arbitrary family
of elements of $\Sigma$, then also
$$
\coprod_{i\in I}\phi_i:\coprod_{i\in I}A_i\to\coprod_{i\in I}B_i
$$
lies in $\Sigma$.
\item
Any retract of any element of $\Sigma$ lies also in
$\Sigma$ (see definition \ref{def_fibrant-Kan}(iii)).
\end{enumerate}

(ii)\ \ 
The {\em saturation} of $\Sigma$ is the intersection of all
saturated sets of monomorphisms of $s.\Set$ containing $\Sigma$.
\end{definition}

\sset\subsubsection{}\label{subsec_sigma-of-psi}
With this terminology, for any morphism $\psi$ of $s.\Set$,
let $\Sigma(\psi)$ be the set of monomorphisms $\phi$ of
$s.\Set$ such that $\psi$ has the right lifting property
with respect to $\phi$. We point out the following simple
observation :

\begin{lemma}\label{lem_satur-RLP}
With the notation of \eqref{subsec_sigma-of-psi}, the
following holds :
\begin{enumerate}
\item
$\Sigma(\psi)$ is saturated, for any morphism $\psi$
of $s.\Set$.
\item
For any cartesian square of simplicial sets
$$
\xymatrix{ X' \ar[r] \ar[d]_{\psi'} & X \ar[d]^\psi \\
           Y' \ar[r] & Y
}$$
we have $\Sigma(\psi)\subset\Sigma(\psi')$.
\item
If $\psi'$ is any retract of $\psi$, then
$\Sigma(\psi)\subset\Sigma(\psi')$ (see definition
{\em\ref{def_fibrant-Kan}(iii)}).
\item
Especially, in both {\em(ii)} and {\em(iii)}, if $\psi$
is a fibration, the same holds for $\psi'$, and for any
$y\in Y[0]$, the fibre $\psi^{-1}(y)$ is a fibrant
simplicial set (see example {\em\ref{ex_simplicial-sets}(vi)}).
\end{enumerate}
\end{lemma}
\begin{proof} Left to the reader.
\end{proof}

\sset\subsubsection{}
The Kan extension condition isolates a class of simplicial
sets on which the standard constructions of homotopy theory
can be carried out. Many of these constructions come down
to exhibiting certain morphisms from a product of the type
$\bDelta_k\times\bDelta_1$ (for variable $k\in\N$) to given
simplicial sets. To manipulate with ease such basic products,
it will be useful to have at our disposal presentations for
them, in the vein of example \ref{ex_simplicial-sets}(iv,v).

To this aim, notice that $\Delta$ is a full subcategory of
the category $\Poset$ of partially ordered sets (notation
of example \ref{ex_universe}(iii)). We may then extend the
Yoneda embedding of $\Delta$ into $s.\Set$ to a well
defined functor
\set\begin{equation}\label{eq_POSet-into-sSet}
\Poset\to s.\Set
\qquad
(P,\leq)\mapsto\bDelta_P.
\end{equation}
Namely, for every partially ordered set $(P,\leq)$, let
$h_P:\Poset^o\to\Set$ be the image of $P$ under the Yoneda
embedding of the category $\Poset$; then $\bDelta_P$ is
the restriction of $h_P$ to the subcategory $\Delta^o$.
We notice :

\begin{lemma}\label{lem_extend-to-infty}
The functor \eqref{eq_POSet-into-sSet} is fully faithful
and commutes with all limits.
\end{lemma}
\begin{proof} The Yoneda embedding commutes with
representable limits (corollary \ref{cor_pre-misc}(vi)),
and the same holds for the restriction functor
$\Poset^\wedge\to s.\Set$, since limits are computed
argumentwise in both of these categories (corollary
\ref{cor_pre-misc}(ii)). To see that \eqref{eq_POSet-into-sSet}
is faithful, it suffices to remark that there is a natural
isomorphism
\set\begin{equation}\label{eq_natural-poset}
\bDelta_P[0]\isom P
\qquad
\text{for every partially ordered set $(P,\leq)$}
\end{equation}
which induces a natural identification : $\bDelta_\phi=\phi$,
for every morphism $\phi:(P,\leq)\to(Q,\leq)$ of partially
ordered sets. To check that \eqref{eq_POSet-into-sSet} is
full, let $(P,\leq)$ and $(Q,\leq)$ be any two objects of
$\Poset$, and $f:\bDelta_P\to\bDelta_Q$ a morphism of
simplicial sets. The natural identification
\eqref{eq_natural-poset} yields a map $\phi:=f[0]:P\to Q$,
and it suffices to check that $\phi$ is a map of ordered
sets such that $\bDelta_\phi=f$. Thus, let $x_0,x_1\in P$ be
any two elements with $x<y$; there follows a unique morphism
$\psi:\bDelta_1\to P$ of partially ordered sets such that
$\psi(i)=x_i$ for $i=0,1$. Then $\psi\in\bDelta_P[1]$,
and we may set
$$
\psi':=f[1](\psi)\in\bDelta_Q[1]
\qquad
y_i:=\bDelta_Q[\eps_{i-1}](\psi')
\qquad
\text{for $i=0,1$}
$$
(notation of \eqref{subsec_face-and-degeracies}). Tracing
back the definitions, we see that $\phi(x_i)=y_i$ for $i=0,1$,
and the existence of $\psi'$ tells us that $y_0\leq y_1$,
{\em i.e.} $\phi$ is a morphism in $\Poset$, as required.
It remains only to show that $\bDelta_\phi[n]=f[n]$ for
every $n\in\N$. However, let $x:[n]\to P$ be any element
of $\bDelta_P[n]$ and set $y:=f[n](x)$, so the image of
$x$ (resp. of $y$) is a totally ordered sequence
$(x_0,\dots,x_n)$ of elements of $P$ (resp. $(y_0,\dots,y_n)$
of $Q$); we have to check that $\phi(x_k)=y_k$ for every
$k=0,\dots,n$. Now, for any such $k\leq n$, let
$\beta_k:[0]\to[n]$ be the unique map such that $\beta_k(0)=k$;
the correspondence \eqref{eq_natural-poset} identifies
$x_k$ with $\bDelta_P[\beta_k](x)$ and $y_k$ with
$\bDelta_Q[\beta_k](y)$, whence the contention (details
left to the reader).
\end{proof}

\sset\subsubsection{}\label{subsec_present-product}
The category of partially ordered sets is complete;
especially, all finite products
$$
[k_0]\times\cdots\times[k_n]
$$
are representable in $\Poset$ : explicitly, the partial
ordering on such a product is defined by declaring that
$$
\underline a:=(a_0,\dots,a_n)\leq\underline b:=(b_0,\dots,b_n)
\qquad
\text{if and only if $a_i\leq b_i$ for every $i=0,\dots,n$}
$$
for every pair of elements
$\underline a,\underline b\in[k_0]\times\cdots\times[k_n]$.
Now, for any $k,r\in\N$ with $r\leq k$, consider the morphism
in $\Poset$
$$
\phi_{k,r}:[k+1]\to[k]\times[1]
\quad
\text{such that}
\quad
\phi_{k,r}(i)=\left\{\begin{array}{lll}
(0,i) & \quad & \text{for $i=0,\dots,r$} \\
(1,i-1) &        & \text{for $i=r+1,\dots,k+1$}.
                    \end{array}\right.
$$
It is easily seen that every morphism $[t]\to[k]\times[1]$
in $\Poset$ (for any $t\in\N$) factors through some
$\phi_{k,r}$, so we get an epimorphism of simplicial sets
$$
\phi_k:[k]\times\bDelta_{k+1}\to\bDelta_k\times\bDelta_1
$$
whose restriction to each subobject $\{r\}\times\bDelta_{k+1}$
is the morphism $\bDelta_{\phi_{k,r}}$ (notation of
\eqref{eq_POSet-into-sSet}).
Now, let $r,s$ be any two integers such that $0\leq r<s\leq k$;
obviously, the intersection of the images of $\phi_{k,r}$ and
$\phi_{k,s}$ is the subset
$$
\{(0,0),\dots,(0,r),(1,s),\dots,(1,k)\}
$$
so we get a cartesian diagram in $\Poset$
$$
{\diagram [t] \ar[r]^-{\eps_{t,r+1}^{s-r}} \ar[d]_{\eps_{t,r+1}^{s-r}}
        & [k+1] \ar[d]^{\phi_{k,r}} \\
          [k+1] \ar[r]^-{\phi_{k,s}} & [k]\times[1]
\enddiagram}
\qquad
\text{with $t:=k+1-s+r$}
$$
(notation of example \eqref{ex_face-and-deg}(ii)). Since
\eqref{eq_POSet-into-sSet} commutes with fibre products,
and $\phi_k$ is an effective epimorphism (by
\eqref{subsec_eff-sset-epis}), we conclude that the diagram
\set\begin{equation}\label{eq_old-present-arm}
\displaystyle\bigcup_{0\leq r<s\leq k}\{(r,s)\}\times\bDelta_{k+1-s+r}
\xymatrix{\ar@<.5ex>[r]^-{\eps'_\bullet} \ar@<-.5ex>[r]_-{\eps''_\bullet} &
[k]\times\bDelta_{k+1} \ar[r]^-{\phi_k} & \bDelta_k\times\bDelta_1
}\end{equation}
identifies $\bDelta_k\times\bDelta_1$ with the coequalizer of
$\eps'_\bullet$ and $\eps''_\bullet$, where $\eps'_\bullet$
(resp. $\eps''_\bullet$) is the morphism whose restrictions
to each subobject $\{(r,s)\}\times\bDelta_t$ agrees with
$\eps^{s-r}_{t,r+1}:\{(r,s)\}\times\bDelta_t\to\{r\}\times\bDelta_{k+1}$
(resp. with
$\eps^{s-r}_{t,r+1}:\{(r,s)\}\times\bDelta_t\to\{s\}\times\bDelta_{k+1}$).

However, \eqref{eq_old-present-arm} can be simplified as follows.
Notice that if $i_1,\dots,i_r$ is any sequence of non-negative
integers such that $1\geq i_{n+1}-i_n\geq 0$ for every
$n=0,\dots,r-1$, then
$$
\eps_{i_r}\circ\cdots\circ\eps_{i_1}=\eps_{t,i_1}^r
\qquad
\text{for every $t\in\N$}.
$$
From this observation, it is easily seen that the terms
$\{(r,s)\}\times\Delta_{k+1-s-r}$ in \eqref{eq_old-present-arm}
with $s>r+1$ are redundant; so we arrive at the presentation
\set\begin{equation}\label{eq_present-arm}
[k-1]\times\bDelta_k
\xymatrix{\ar@<.5ex>[r]^-{\eps'_\bullet} \ar@<-.5ex>[r]_-{\eps''_\bullet} &
[k]\times\bDelta_{k+1} \ar[r]^-{\phi_k} & \bDelta_k\times\bDelta_1
}\end{equation}
where each term $\{r\}\times\bDelta_k$ corresponds to the
term $\{(r,r+1)\}\times\bDelta_k$ of \eqref{eq_old-present-arm},
and the morphisms $\eps'_\bullet$, $\eps''_\bullet$ are redefined
accordingly.

\sset\subsubsection{}\label{subsec_saturate}
We consider now the following sets of monomorphisms of $s.\Set$ :
\begin{itemize}
\item
$\Sigma_1$ is the set of all morphisms
$\iota_k^n:\bLambda_k^n\to\bDelta_k$, for every $k,n\in\N$ with
$k\geq n$.
\item
$\Sigma_2$ is the set of all morphisms
$$
i^j_k:
(\partial\bDelta_k\times\bDelta_1)\cup(\bDelta_k\times\bLambda^j_1)
\to\bDelta_k\times\bDelta_1
\qquad
\text{for every $k\in\N$ and $j=0,1$}.
$$
\item
$\Sigma_3$ is the set of all monomorphisms of the form
$$
i^j_{K,L}:
(K\times\bDelta_1)\cup(L\times\bLambda_1^j)\to L\times\bDelta_1
$$
where $K\to L$ is an arbitrary monomorphism of $s.\Set$, and
$j=0,1$.
\end{itemize}

With this notation, we may now state :

\begin{proposition}\label{prop_saturate}
The sets $\Sigma_1$, $\Sigma_2$ and $\Sigma_3$ have the
same saturation.
\end{proposition}
\begin{proof} In order to check that the saturation of
$\Sigma_1$ contains $\Sigma_2$, fix $k\in\N$; if $j=0$
(resp. if $j=1$), for every $r=0,\dots,k+1$, denote by
$A_r\subset\bDelta_k\times\bDelta_1$  the image of the
restriction of $\phi_k$ to the subobject
$\{0,\dots,r\}\times\bDelta_{k+1}$ (resp. the subobject
$\{k-r+1,\dots,k+1\}\times\bDelta_{k+1}$), and define $B_r$
as the fibre product in the cartesian diagram 
$$
\xymatrix{ B_r \ar[r] \ar[d] &
(\partial\bDelta_k\times\bDelta_1)\cup(\bDelta_k\times\bLambda^j_1)
\ar[d]^{i^j_k} \\
A_r \ar[r] & \bDelta_k\times\bDelta_1.
}$$
Set also $A_{-1}:=s.\emptyset$ (the initial object of $s.\Set$).
Now, a straightforward induction reduces to checking that the
natural morphisms
$$
B_k\amalg_{B_r}A_r\to B_k\amalg_{B_{r+1}}A_{r+1}
\qquad
\text{for $r=-1,\dots,k-1$}
$$
lie in the saturation of $\Sigma_1$. To this aim, it
suffices to show that the same holds for the natural
morphisms
$$
\tau_r:B_{r+1}\amalg_{B_r}A_r\to A_{r+1}
\qquad
\text{for $r=-1,\dots,k-1$}.
$$
This, in turns, follows immediately from :

\begin{claim}
(i)\ \
The natural isomorphism $\{0\}\times\bDelta_{k+1}\isom A_0$
identifies $B_0$ with $\{0\}\times\bLambda_{k+1}^1$.

(ii)\ \
For every $r=0,\dots,k-1$ there is a commutative diagram
in $s.\Set$
$$
\xymatrix{
\bDelta_k \ar[r]^-{\eps_{r+1}} \ar[d]_{\eps_{r+1}} &
\bLambda^{r+2}_{k+1} \ar[r]^-{\iota^{r+2}_{k+1}} \ar[d] &
\bDelta_{k+1} \ar[d] \\
A_r \ar[r] & B_{r+1}\amalg_{B_r}A_r \ar[r]^-{\tau_r} & A_{r+1}
}$$
whose two square subdiagrams are cocartesian.
\end{claim}
\begin{pfclaim} To ease notation, we shall identify
$\{r\}\times\bDelta_k$ with its image under $\phi_k$,
which is a subobject of $\bDelta_k\times\bDelta_1$, for
every $r=0,\dots,k$. Notice that $B_{r+1}\amalg_{B_r}A_r$ is
the smallest subobject of $\bDelta_k\times\bDelta_1$ that
contains $A_r$ and $Z_{r+1}:=(\{r+1\}\times\bDelta_{k+1})\cap B_k$.
However, a morphism $[t]\to B_k$ is the same as an increasing
sequence of elements of $[k]\times[1]$
\set\begin{equation}\label{eq_a-sequence}
(a_n,b_n)
\qquad
n=0,\dots,t
\end{equation}
such that :
\begin{enumerate}
\alphaenu
\item
either the cardinality of the set $\{a_0,\dots,a_t\}$
is strictly less than $k+1$
\item
or else $b_n\neq j$ for every $n=0,\dots,t$.
\end{enumerate}
On the other hand, recall that the inclusion
$\{r+1\}\times\bDelta_{k+1}\subset\bDelta_k\times\bDelta_1$
corresponds to the inclusion map of partially ordered sets
$\phi_{k,r+1}$ (notation of \eqref{subsec_present-product}).
It follows that an injective morphism $\bDelta_t\to Z_{r+1}$
corresponds to a strictly increasing sequence
\eqref{eq_a-sequence} contained in
$$
\{(0,0),\dots,(0,r+1),(1,r+1),\dots,(1,k)\}
$$
and fulfilling either of the foregoing conditions (a) or (b).
Taking $r=-1$, we already see that (i) holds.
If $r\geq 0$, take $t:=k$, and notice that no such sequence can
have $a_r\neq r+1$, so we see that of all the $k+2$ monomorphisms
$$
\one_{\{r+1\}}\times\eps_j:
\{r+1\}\times\bDelta_k\to\{r+1\}\times\bDelta_{k+1}
$$
only those with $j\neq r+1,r+2$ factor through $Z_{r+1}$.
On the other hand, from the presentation \eqref{eq_present-arm}
we see that the intersection of $A_r$ and
$\{r+1\}\times\bDelta_{k+1}$ is the image of
$\{r\}\times\bDelta_k$, which maps to both of them via
the morphism $\eps_{r+1}$. From this, we conclude already that
the left square subdiagram of (ii) is indeed cocartesian,
if we let the central vertical arrow to be the morphism that
naturally identifies $\bLambda^{r+2}_{k+1}$ with the subobject
$\{r+1\}\times\bLambda^{r+2}_{k+1}$ of $\{r+1\}\times\bDelta_{k+1}$.
By the same token, it is clear the right square subdiagram is
likewise cocartesian, provided we let the right-most vertical
arrow to be the morphism that naturally identifies $\bDelta_{k+1}$
with $\{r+1\}\times\bDelta_{k+1}$.
\end{pfclaim}

Next, let us check that the saturation of $\Sigma_2$
contains $\Sigma_3$. Indeed, let $\mu:K\to L$ be any
monomorphism of $s.\Set$, and fix $j\in\{0,1\}$;
according to example \ref{ex_simplicial-sets}(iii),
there exist :
\begin{itemize}
\item
A countable system of monomorphisms
$$
K_{-1}:=K\xrightarrow{\ \mu_{-1}\ }K_0\to
\cdots\to K_n\xrightarrow{\ \mu_n\ }K_{n+1}\to\cdots
$$
such that $\mu$ is isomorphic to the induced morphism
$$
K\to\colim_{n\in\N}K_n.
$$
\item
For every $n\in\N$, a set $I_n$ and an epimorphism
$K_{n-1}\amalg(I_n\times\bDelta_n)\to K_n$
whose restriction to $K_{n-1}$ agrees with $\mu_{n-1}$.
\end{itemize}
There follows, for every integer $n\in\N$, a cocartesian diagram
$$
\xymatrix{
(K_{n-1}\times\bDelta_1)\cup(K_n\times\bLambda^j_1)
\ar[r] \ar[d]_-{i^j_{K_{n-1},K_n}} &
(K_{n-1}\times\bDelta_1)\cup(L\times\bLambda^j_1) \ar[d]^{\phi_n} \\
K_n\times\bDelta_1 \ar[r] &
(K_n\times\bDelta_1)\cup(L\times\bLambda^j_1)
}$$
and the system $(\phi_n~|~n\in\N)$ induces a morphism
$$
(K\times\bDelta_1)\cup(L\times\bLambda^j_1)\to
\colim_{n\in\N}(K_n\times\bDelta_1)\cup(L\times\bLambda^j_1)
=L\times\bDelta_1
$$
which is isomorphic to $i^j_{K,L}$.
Thus, we may assume that there exist a set $I$ and an
epimorphism $K\amalg(I\times\bDelta_n)\to L$ for some
$n\in\N$, and we shall argue by induction on $n$.
Set $Q:=(I\times\bDelta_n)\times_LK$ and
notice that $Q$ is a subobject of $I\times\bDelta_n$.
Also, since all epimorphisms are effective, the induced
commutative diagram
$$
\xymatrix{ Q \ar[r] \ar[d] & K \ar[d] \\
I\times\bDelta_n \ar[r] & L
}$$
is cocartesian. Furthermore, since all colimits are
universal in $s.\Set$, we deduce that
$$
Q=\coprod_{i\in I_n}Q_i
\qquad\text{where}\qquad
Q_i:=(\{i\}\times\bDelta_n)\cap Q
\quad
\text{for every $i\in I$}.
$$
We are then reduced to the case where $L=\bDelta_n$.
If $n=0$, then clearly $K$ is either $\bDelta_0$
or $\bDelta_{-1}$, and in either case we have
$i^j_{K,L}\in\Sigma_2$.

Next, suppose that $n>0$, and that the assertion is already
known for every integer $k<n$. If $K=\bDelta_n$, we are
done; otherwise, it is easily seen that the inclusion
$K\to\bDelta_n$ factors through $\partial\bDelta_n$, and
arguing as in the foregoing, we are reduced to considering
the morphisms $i^j_{K,\partial\bDelta_n}$ and
$i^j_{\partial\bDelta_n,\bDelta_n}$. The latter lies in $\Sigma_2$,
and from example \ref{ex_simplicial-sets}(iv) it is clear
that there is an epimorphism
$K\amalg([n]\times\bDelta_{n-1})\to\partial\bDelta_n$, so
$i^j_{K,\partial\bDelta_n}$ lies in the saturation of $\Sigma_2$,
by inductive assumption.

Lastly, we show that the saturation of $\Sigma_3$ contains
$\Sigma_1$. To this aim, it suffices to construct a
commutative diagram
$$
\xymatrix{ \bLambda^r_k \ar[r] \ar[d]_{\iota^r_k} &
(\bLambda^r_k\times\bDelta_1)\cup(\bDelta_k\times\bLambda_1^j)
\ar[r] \ar[d]^{i^j_{\bLambda^r_k,\bDelta_k}} & \bLambda^r_k \ar[d]^{\iota_k^r} \\
\bDelta_k \ar[r]^-i & \bDelta_k\times\bDelta_1 \ar[r]^-p & \bDelta_k
}$$
for every $n,r\in\N$ with $r\leq k$, such that
$p\circ i=\one_{\bDelta_k}$. However, notice that a morphism
$[t]\to(\bLambda^r_k\times\bDelta_1)\cup(\bDelta_k\times\bLambda_1^j)$
(for any $t\in\N$) is the same as an increasing sequence
\eqref{eq_a-sequence} of elements of $[k]\times[1]$ such that
\begin{enumerate}
\alphaenu
\item
either $b_n\neq j$ for every $n=0,\dots,t$
\item
or else, the cardinality of of the set $\{a_0,\dots,a_t\}\cup\{r\}$
is strictly less than $k+1$.
\end{enumerate}
Now, any such $p$ and $i$ shall be the images, under the functor
\eqref{eq_POSet-into-sSet}, of corresponding maps of ordered sets
$$
[k]\xrightarrow{\ i^*\ }[k]\times[1]\xrightarrow{\ p^*\ }[k]
\qquad
\text{such that $p^*\circ i^*=\one_{[k]}$}.
$$
We obtain a suitable $j\in\{0,1\}$ and suitable $p^*$, $i^*$
by the following rule. If $r<k$, we let $j:=1$, and define
$p^*$ and $i^*$ as the maps such that
$$
i^*(a)=(a,1)
\quad
\text{for every $a=0,\dots,k$}
\qquad
p^*(a,b)=\left\{\begin{array}{lll}
              a & \quad & \text{if $b=1$ or $a\leq r$} \\
              r &       & \text{otherwise}.
              \end{array}\right.
$$
If $r=k$, we let $j:=0$, and set
$$
i^*(a):=(a,0)
\quad
\text{for every $a=0,\dots,k$}
\qquad
p^*(a,b)=\left\{\begin{array}{lll}
              a & \quad & \text{if $b=0$} \\
              r &       & \text{otherwise}.
              \end{array}\right.
$$
The reader may easily check that the resulting maps $i$ and
$p$ will do.
\end{proof}

\begin{definition} A monomorphism in $s.\Set$ is an
{\em anodyne extension}, if it lies in the saturation
of the set $\Sigma_1$.
\end{definition}

In view of lemma \ref{lem_satur-RLP}(i), we see that a
fibration has the right lifting property with respect
to every anodyne extension. We also notice :

\begin{corollary}\label{cor_saturate}
Let $K\to L$ be any anodyne extension, and $A\to B$ any
monomorphism. Then the induced monomorphism
$$
(K\times B)\cup(L\times A)\to L\times B
$$
is an anodyne extension.
\end{corollary}
\begin{proof} Let $\Sigma$ be the set of monomorphisms
$K'\to L'$ such that the induced morphism
$$
(K'\times B)\cup(L'\times A)\to L'\times B
$$
is an anodyne extension. It is easily seen that $\Sigma$ is
saturated (details left to the reader; cp. the proof of
proposition \ref{prop_saturate}); taking into account
proposition \ref{prop_saturate}, it then suffices to check
that $\Sigma$ contains the morphisms $i^j_{X,Y}$, where
$f:X\to Y$ is an arbitrary monomorphism (notation of
\eqref{subsec_saturate}). However, for any such $f$ we
have a commutative diagram
$$
\xymatrix{
(((X\times\bDelta_1)\cup(Y\times\bLambda_1^j))\times B)
\cup((Y\times\bDelta_1)\times A) \ar[rrr] \ar[d] & & &
((Y\times\bDelta_1)\times B \ar[d] \\
(((X\times B)\cup(Y\times A))\times\bDelta_1)
\cup((Y\times B)\times\bLambda_1^j)
\ar[rrr]^-{i^j_{(X\times B)\cup(Y\times A),Y\times B}} & & &
(Y\times B)\times\bDelta_1
}$$
whose vertical arrows are isomorphisms. Then the contention
follows by appealing again to proposition \ref{prop_saturate}.
\end{proof}

\begin{theorem}\label{th_E-infinity}
For every morphism $f:X\to Y$ of simplicial sets there
exists a commutative diagram
$$
{\diagram X \ar[rr]^{i_f} \ar[rd]_f & & E_f \ar[ld]^{p_f} \\
& Y
\enddiagram}
\qquad
\text{in $s.\Set$}
$$
such that $i_f$ is an anodyne extension, and $p_f$ is a
fibration.
\end{theorem}
\begin{proof} For every $k,n\in\N$ with $k\geq\min(1,n)$,
consider the set $\cL^n_k$ of all commutative diagrams
$$
\xymatrix{\bLambda^n_k \ar[r] \ar[d]_{\iota^n_k} &
X \ar[d]^f \\
\bDelta_k \ar[r] & Y
}$$
and set
$$
\bLambda_f:=
\coprod_{n\in\N}\ \coprod_{k\geq\min(1,n)}\!\!\cL^n_k\times\bLambda^n_k
\qquad
\bL_f:=
\coprod_{n\in\N}\ \coprod_{k\geq\min(1,n)}\!\!\cL^n_k\times\bDelta_k.
$$
There follows a natural morphism
$$
\xymatrix{ \bLambda_f \ar[r]^-g \ar[d]_\iota & X \ar[d]^f \\
\bL_f \ar[r] & Y
}$$
with $\iota$ an anodyne extension. Define $E^0_f$ as the
push-out in the induced cocartesian diagram
$$
\xymatrix{
\bLambda_f \ar[r]^-g \ar[d]_\iota & X \ar[d]^{i^0_f} \\
\bL_f \ar[r] & E^0_f
}$$
so that $i^0_f$ is anodyne as well, and we have a commutative
diagram
$$
\xymatrix{
X \ar[rr]^-{i^0_{(X,f)}} \ar[rd]_f & & E^0_f \ar[ld]^{p^0_f} \\
& Y.
}$$
Clearly, the rule $(X,f)\mapsto(E^0_f,p^0_f)$ extends to a
well defined functor
$$
\underline E^0:s.\Set/Y\to s.\Set/Y
$$
and the rule $(X,f)\mapsto i^0_{(X,f)}$ yields a natural
transformation
$$
i^0:\one_{s.\Set/Y}\Rightarrow\underline E^0.
$$
Thus, we may define inductively :
$$
\underline E^n(X,f):=\underline E^0(\underline E^{n-1}(X,f))
\qquad
\text{for every integer $n\geq 1$}
$$
so $\underline E^n(X,f)$ is a pair $(E^n_f,p^n_f)$, and
we get as well as a natural transformation
$$
i^n_{(X,f)}:=i^0_{\underline E^{n-1}(X,f)}:E^{n-1}_f\to E^n_f
\qquad
\text{for every integer $n\geq 1$}
$$
which is again an anodyne extension. We set
$$
E_f:=\colim_{n\in\N} E^n_f
$$
where the transition maps in the colimit are given by
the system of morphisms $(i^n_f~|~n\in\N)$. The colimit
of the system of morphism $(p^n_f~|~n\in\N)$ is then a
morphism $p_f:E_f\to Y$; moreover, the induced morphism
$i_f:X\to E_f$ is anodyne, and obviously $p_f\circ i_f=f$.
To conclude the proof, it suffices to show that $p_f$
is a fibration. Thus, consider any commutative diagram
$$
\xymatrix{
\bLambda^n_k \ar[r]^-g \ar[d]_{\iota^n_k} &
E_f \ar[d]^{p_f} \\
\bDelta_k \ar[r]^-h & Y.
}$$
It follows easily from example \ref{ex_simplicial-sets}(v)
that -- for $r\in\N$ large enough -- $g$ factors through a
morphism $g_r:\bLambda^n_k\to E^r_f$ and the natural
morphism $E^r_f\to E_f$, whence a commutative diagram
$$
\xymatrix{
\bLambda^n_k \ar[r]^-{g_r} \ar[d]_{\iota^n_k} &
E^r_f \ar[d]^{p^r_f} \ar[r]^-{i_f^{r+1}} & E^{r+1}_f \ar[d]^{p^{r+1}_f} \\
\bDelta_k \ar[r]^-h & Y \rdouble & Y.
}$$
But then, by construction of $E^{n+1}_f$, we see that
$h$ lifts to a morphism $h':\bDelta_k\to E^{n+1}_f$ such
that $p^{r+1}_f\circ h'=h$ and
$h'\circ\iota^n_k=i_f^{r+1}\circ g_r$. The composition
of $h'$ with the natural morphism $E^{r+1}_f\to E_f$ is
a morphism $h'':\bDelta_k\to E_f$ such that $p_f\circ h''=h$
and $h''\circ\iota^n_k=g$, as required.
\end{proof}

A useful way to produce new fibrant simplicial sets out of old
ones is provided by the general construction introduced in the
following :

\begin{definition}
Let $A$ and $B$ be any two simplicial sets; the
{\em function complex}
$$
s.\cHom(A,B)
$$
associated with $A$ and $B$ is the simplicial set given by the rule :
$$
n\mapsto\Hom_{s.\Set}(\bDelta_n\times A,B)
\qquad
\phi\mapsto\Hom_{s.\Set}(\bDelta_\phi\times\one_A,B)
$$
for every $n\in\N$ and every morphism $\phi$ of $\Delta$.
\end{definition}

\begin{remark}\label{rem_simpl-hom}
(i)\ \
Notice that the system of isomorphisms
$A[n]\isom\Hom_{s.\Set}(\bDelta_n,A)$ given by Yoneda's lemma
for every $n\in\N$, amounts to a natural identification
$$
\iota_A:A\isom s.\cHom(\bDelta_0,A)
\qquad
\text{for every $A\in\Ob(s.\Set)$}.
$$

(ii)\ \
Also, if $A$, $B$ and $C$ are any three simplicial sets,
there is a natural transformation
$$
\tau_{A,B,C}:s.\cHom(A,B)\to s.\cHom(A\times C,B\times C)
$$
which, to any $n\in\N$, assigns the map
$$
\Hom_{s.\Set}(\bDelta_n\times A,B)\to
\Hom_{s.\Set}(\bDelta_n\times A\times C,B\times C)
\qquad
\phi\mapsto\phi\times\one_C.
$$

(iii)\ \
There is a natural {\em evaluation morphism}
$$
\ev_{A,B}:A\times s.\cHom(A,B)\to B
\qquad
\text{for every $A,B\in\Ob(s.\Set)$}
$$
which, to every $n\in\N$, assigns the map
$$
A[n]\times\Hom_{s.\Set}(\bDelta_n\times A,B)\to B[n]
\qquad
(a,f)\mapsto f[n](\one_{\bDelta_n},a).
$$
Indeed, for any morphism $\phi:[k]\to[n]$ in $\bDelta$,
and $(a,f)$ any element of $(A\times\cHom(A,B))[n]$ we
may compute
$$
\begin{aligned}
B[\phi]\circ\ev_{A,B}[n](a,f)
= &\, B[\phi]\circ f[n](\one_{\bDelta_n},a) \\
= &\, f[k]\circ(\bDelta_n[\phi]\times A[\phi])(\one_{\bDelta_n},a) \\
= &\, f[k](\phi,A[\phi](a)) \\
= &\, f[k]\circ(\bDelta_\phi\times\one_A)[k](\one_{\bDelta_k},A[\phi](a)) \\
= &\, (f\circ(\bDelta_\phi\times\one_A))[k](\one_{\bDelta_k},A[\phi](a)) \\
= &\, \ev_{A,B}[k](A[\phi](a),f\circ(\bDelta_\phi\times\one_A)) \\
= &\, \ev_{A,B}[k]\circ(A\times\cHom(A,B))[\phi](a,f)
\end{aligned}
$$
which shows that $\ev_{A,B}$ is a morphism of simplicial sets.

(iv)\ \
Clearly $s.\Set$ is a tensor category, with tensor product
given by the product of simplicial sets, and with $\bDelta_0$
as unit object. We claim that the functor
$$
s.\cHom:(s.\Set)^o\times s.\Set\to s.\Set
$$
is an internal $\Hom$ functor for $(s.\Set,\times,\bDelta_0)$.
Indeed, given $A,B,C\in\Ob(s.\Set)$ and any morphism
$f:A\times B\to C$, we obtain a morphism
$$
A\xrightarrow{\ \iota_A\ }s.\cHom(\bDelta_0,A)
\xrightarrow{\ \tau_{A,B,C}\ }s.\cHom(B,A\times B)
\xrightarrow{\ \cHom(B,f)\ }s.\cHom(B,C).
$$
Conversely, given a morphism $g:A\to\cHom(B,C)$, we get a
morphism
$$
A\times B\xrightarrow{\ g\times\one_B\ }
s.\cHom(B,C)\times B\isom B\times s.\cHom(B,C)
\xrightarrow{\ \ev_{B,C}\ }C.
$$
It is easily seen that these two rules yield mutually
inverse natural bijections, as required (details left
to the reader).
\end{remark}

\sset\subsubsection{}\label{subsec_fibrant-hom}
Let $f:A\to B$ and $g:C\to D$ be morphisms of simplicial
sets. We have a commutative diagram
$$
\xymatrix{ s.\cHom(D,A) \ar[rr]^-{s.\cHom(D,f)}
\ar[d]_{s.\cHom(g,A)} & & s.\cHom(D,B) \ar[d]^{s.\cHom(g,B)} \\
s.\cHom(C,A) \ar[rr]^-{s.\cHom(C,f)} & & s.\cHom(C,B)
}$$
whence an induced morphism in $s.\Set$
\set\begin{equation}\label{eq_fibrant-hom}
s.\cHom(D,A)\to s.\cHom(D,B)\times_{s.\cHom(C,B)}s.\cHom(C,A).
\end{equation}

\begin{proposition}\label{prop_fibrant-hom}
In the situation of \eqref{subsec_fibrant-hom}, suppose
that $f$ is a fibration and $g$ a monomorphism. Then
\eqref{eq_fibrant-hom} is a fibration as well.
\end{proposition}
\begin{proof} The datum of a commutative diagram
$$
\xymatrix{
\bLambda^i_k \ar[r] \ar[d]_{\iota_k^i} & s.\cHom(D,A) \ar[d] \\
\bDelta_k \ar[r] & s.\cHom(D,B)\times_{s.\cHom(C,B)}s.\cHom(C,A)
}$$
is equivalent, by remark \ref{rem_simpl-hom}(iv), to that
of a commutative diagram
$$
\xymatrix{ (\bLambda^i_k\times D)\cup(\bDelta_k\times C)
\ar[r] \ar[d]_g & A \ar[d]^f \\
\bDelta_k\times D \ar[r] & B.
}$$
On the other hand, $f$ has the right lifting property with
respect to $g$ (corollary \ref{cor_saturate}); it follows
easily that \eqref{eq_fibrant-hom} has the right lifting
property with respect to $\iota^i_k$, as stated.
\end{proof}

\begin{corollary}\label{cor_fibrant-hom}
Let $A$ be any fibrant simplicial set. We have :
\begin{enumerate}
\item
$s.\cHom(D,A)$ is fibrant, for any simplicial set $D$.
\item
Any monomorphism $D\to D'$ of simplicial sets induces a fibrant
morphism
$$
s\cHom(D',A)\to s.\cHom(D,A).
$$
\end{enumerate} 
\end{corollary}
\begin{proof}(i): Take $B:=\bDelta_0$, $C:=\bDelta_{-1}$,
and let $f:A\to B$, $g:C\to D$ be the unique morphisms
of simplicial sets. In view of proposition \ref{prop_fibrant-hom},
it suffices to notice the natural identifications
$$
s.\cHom(D,\bDelta_0)\isom\bDelta_0\isom s.\cHom(\bDelta_{-1},D)
$$
for every simplicial set $A$. The proof of (ii) is similar :
details left to the reader.
\end{proof}

\sset\subsubsection{}\label{subsec_homot-rel-to-sub}
Let $f,g:X\to Y$ be two morphisms of simplicial sets.
Recall (see remark \ref{rem_mixed-simpl-tensors}(v))
that a {\em homotopy $u$ from $f$ to $g$} is the datum
of a commutative diagram of simplicial sets :
$$
\xymatrix{ X\times\bDelta_0 \ar[rr]^-{\one_X\times\eps_0}
           \ddouble & & X\times\bDelta_1 \ar[d]^u & & \ddouble 
           X\times\bDelta_0 \ar[ll]_{\one_X\times\eps_1} \\
X \ar[rr]^-f & & Y & & X \ar[ll]_-g.
}$$

\begin{definition}\label{def_homotop_rel}
In the situation of \eqref{subsec_homot-rel-to-sub},
suppose that $i:X'\to X$ is a monomorphism of simplicial
sets, such that $f\circ i=g\circ i$. We say that $u$ is
a {\em homotopy from $f$ to $g$ relative to $X'$}, if
$u$ is a homotopy from $f$ to $g$, and the following
diagram commutes :
$$
\xymatrix{ X'\times\bDelta_1 \ar[d]_{i\times\one_{\bDelta_1}}
           \ar[r] & X' \ar[d]^{f\circ i} \\
           X \times\bDelta_1 \ar[r]^-u & Y
}$$
where the top horizontal arrow is the natural projection.
\end{definition}

\begin{remark}\label{rem_constant-homot}
Let $f,g:X\to Y$ be any two morphisms of simplicial sets.

(i)\ \
Of course, a homotopy from $f$ to $g$ can be viewed as a
homotopy relative to the sub-object $\bDelta_{-1}\subset X$,
so definition \ref{def_homotop_rel} generalizes the previous
notion. Notice also that to $f$ one may attach a
{\em constant homotopy}
$$
X\times\bDelta_1\xrightarrow{\ \one_X\times\eta_0\ }
X\xrightarrow{\ f\ }Y.
$$

(ii)\ \
Let $u$ be any homotopy from $f$ to $g$ relative to
a subobject $X'$ of $X$, and $h:Y\to Z$ any morphism
of simplicial sets. Then clearly $h\circ u$ is a
homotopy from $h\circ f$ to $h\circ g$ relative to $X'$.
Likewise, if $t:Z\to X$ is any other morphism, then
$u\circ(t\times\bDelta_1)$ is a homotopy from $f\circ t$
to $g\circ t$ relative to $Z\times_XX'$.
\end{remark}

\sset\subsubsection{}\label{eqref_homot-relation}
Let $A$ be any simplicial set; on the set
$\Hom_{s.\Set}(\bDelta_0,A)\simeq A[0]$ we consider the binary
relation $\sim$ such that
$$
a\sim b
\qquad\Leftrightarrow\qquad
\text{there exists a homotopy from $a$ to $b$}.
$$

\begin{lemma}\label{lem_homotopy-vertex}
With the notation of \eqref{eqref_homot-relation},
suppose that $A$ is fibrant. Then $\sim$ is an equivalence
relation on $A[0]$.
\end{lemma}
\begin{proof} Clearly $a\sim b$ if and only if there exists
a morphism $u:\Delta_1\to X$ such that $u\circ\eps_0=a$ and
$u\circ\eps_1=b$. Taking $u:=a\circ\sigma_0$, where
$\sigma_0:\bDelta_1\to\bDelta_0$ is the unique morphism,
and $a\in A[0]$ any element, we obtain the reflexivity of
$\sim$. Next, suppose that $a\sim b$ and $b\sim c$ for some
$a,b,c\in A[0]$, and pick $u,v:\bDelta_1\to A$ that give
a homotopy from $a$ to $b$, and respectively from $b$ to
$c$. Then it is easily seen that the datum of $u$ and $v$
determines a unique morphism $t:\bLambda_2^1\to A$ such
that $t\circ\eps_0=u$ and $t\circ\eps_2=v$. By assumption,
$t$ extends to a morphism $s:\bDelta_2\to A$, and then
$s\circ\eps_1$ is a homotopy from $a$ to $c$, whence
the transitivity of $\sim$.

Lastly, let $u:\bDelta_1\to A$ be a homotopy from $a$
to $b$, for some $a,b\in A[0]$; we let $t:\bLambda_2^0\to A$
be the unique morphism such that $t\circ\eps_2=u$ and
$t\circ\eps_1=b\circ\sigma_0$. By assumption $t$ extends
to a morphism $s:\bDelta_2\to A$, and it is easily seen
that $s\circ\eps_0$ is a homotopy from $b$ to $a$, which
shows that $\sim$ is also symmetric.
\end{proof}

\sset\subsubsection{}\label{subsec_homotop-rel}
More generally, if $i:K\to X$ is any monomorphism in $s.\Set$,
and $Y$ is any simplicial set, we consider on $\Hom_{s.\Set}(X,Y)$
the binary relation
$$
f\sim_K g
\qquad\Leftrightarrow\qquad
\text{there exists a homotopy from $f$ to $g$ relative to $K$}.
$$
In case $K=s.\emptyset$, we simply write $f\sim g$ instead
of $f\sim_{s.\emptyset}g$.

\begin{theorem}\label{th_hom-is-equiv}
With the notation of \eqref{subsec_homotop-rel}, suppose
that $Y$ is fibrant. Then $\sim_K$ is an equivalence
relation.
\end{theorem}
\begin{proof} The fibres of the map of sets
\set\begin{equation}\label{eq_attention}
\Hom_{s.\Set}(X,Y)\to\Hom_{s.\Set}(K,Y)
\qquad
f\mapsto f\circ i
\end{equation}
give a partition of $\Hom_{s.\Set}(X,Y)$, and if
$f\sim_Kg$, then $f$ and $g$ lie in the same fibre of
\eqref{eq_attention}. Hence, we may restrict attention
to a single fibre of \eqref{eq_attention}, say the fibre
over the morphism $h:K\to Y$. Now, $h$ corresponds to a
morphism $h^*:\bDelta_0\to s.\cHom(K,Y)$ in $s.\Set$, and
likewise, the datum of a pair of morphisms $f,g:X\to Y$
such that $f\circ i=h=g\circ i$ corresponds to that of a
pair of morphisms $f^*,g^*:\bDelta_0\to s.\cHom(X,Y)$ with
$i^*\circ f^*=h^*=i^*\circ g^*$, where
$$
i^*:s.\cHom(X,Y)\to s.\cHom(K,Y)
$$
is the morphism deduced from $i$. Furthermore, a homotopy
from $f$ to $g$ relative to $K$ is the same as a morphism
$t:\bDelta_1\to s.\cHom(X,Y)$ such that $i^*\circ t$ factors
through $h^*$ and such that
$$
t\circ\eps_0=f^*
\qquad
t\circ\eps_1=g^*.
$$
Set $A:=i^{*-1}(h^*)$ (notation of example
\ref{ex_simplicial-sets}(vi)). We conclude that $t$
corresponds to a homotopy between the elements of
$A[0]$ induced by $f^*$ and $g^*$. However, $A$ is
fibrant by corollary \ref{cor_fibrant-hom}(ii) and
lemma \ref{lem_satur-RLP}(iv), so the assertion
follows from lemma \ref{lem_homotopy-vertex}.
\end{proof}

\sset\subsubsection{}
Let $A$ be any simplicial set, $\xi\in A[0]$ any element,
and denote also by $\xi:\bDelta_0\to A$ the corresponding
morphism of simplicial sets.
We shall frequently abuse notation, and denote indifferently
by $\xi$ the unique element of $A[n]$ that lies in the image
of the map $\xi[n]$. We call the pair $(A,\xi)$ a
{\em pointed simplicial set}. With this terminology, we have : 

\begin{definition}\label{def_hom-groups}
Let $(A,\xi)$ be a pointed simplicial set, such that $A$
is fibrant.
\begin{enumerate}
\item
For every integer $n\in\N$, let $A(\xi,n)$ be the set of all
$a\in A[n]$ such that $\partial a$ factors through $\xi$
(notation of example \ref{ex_simplicial-sets}(iv)). We set
$$
\pi_n(A,\xi):=A(\xi,n)/\!\!\sim_{\partial\bDelta_n}
$$
(where $\sim_{\partial\bDelta_n}$ is the homotopy equivalence
relation defined in \eqref{subsec_homotop-rel}), and for every
$n>0$ we call $\pi_n(A,\xi)$ the {\em $n$-th homotopy group}
of $A$.
\item
Notice that $\pi_0(A,\xi)$ is actually independent of
$\xi$, so we shall usually just denote it $\pi_0(A)$.
We say that $A$ is {\em connected} if the cardinality
of $\pi_0(A)$ equals one.
\end{enumerate}
\end{definition}

\begin{lemma}\label{lem_hom-groups}
In the situation of definition {\em\ref{def_hom-groups}},
let $n\in\N$ be any integer and $a\in A(\xi,n)$ any simplex.
The following conditions are equivalent :
\begin{enumerate}
\alphaenu
\item
$a\sim_{\partial\bDelta_n}\xi$.
\item
There exists $b\in A[n+1]$ such that $\partial_ib=\xi$
for every $i=1,\dots,n+1$ and $\partial_0b=a$.
\end{enumerate}
\end{lemma}
\begin{proof} (a)$\Rightarrow$(b): Let
$h:\bDelta_n\times\bDelta_1\to A$ be any morphism in
$s.\Set$. From the presentation \eqref{eq_present-arm}
we see that $h$ is the same as a system $x_0,\dots,x_n$
of elements of $A[n+1]$ such that
$$
\partial_ix_i=\partial_ix_{i-1}
\qquad
\text{for $i=1,\dots,n$}.
$$
Namely, $x_i:=h\circ\bDelta_{\phi_{n,i}}$ for every $i=0,\dots,n$,
where $\phi_{n,i}:[n+1]\to[n]\times[1]$ is defined as in
\eqref{subsec_present-product}. Moreover, recall that
the fully faithful functor \eqref{eq_POSet-into-sSet}
associates with any map of partially ordered sets
$f:[n]\to[n]\times[1]$ a morphism
$\bDelta_f:\bDelta_n\to\bDelta_n\times\bDelta_1$, and
it is easily seen that $\bDelta_f$ factors through the
subobject $(\partial\bDelta_n\times\bDelta_1)\cup
(\bDelta_n\times\partial\bDelta_1)$
if and only if $f$ fulfills the following condition.
Set $f(j):=(r_j,s_j)$ for every $j\in[n]$; then, for
every $j=0,\dots,n-1$ we have either $r_j=r_{j+1}$ or
$s_j=s_{j+1}$ (details left to the reader). It follows
easily that $h$ is a homotopy from $\xi$ to $a$ relative
to $\partial\bDelta_n$ if and only if we have
$$
\partial_0x_0=a
\qquad
\partial_{n+1}x_n=\xi
\qquad\text{and}\qquad
\partial_jx_i=\xi
\quad
\text{for $0<i\leq n$ and $j\neq i,i+1$}.
$$
Consider then the system $y_\bullet:=(y_1,\dots,y_{n+2})$
of elements of $A[n+1]$ such that $y_{n+2}:=\xi$ and $y_i=x_{i-1}$
for every $i=1,\dots,n+1$. A direct calculation shows that
$$
\partial_iy_j=\partial_{j-1}y_i
\qquad
\text{whenever $1\leq i<j\leq n+2$}
$$
so $y_\bullet$ corresponds to a unique morphism
$\phi:\bLambda_{n+2}^0\to A$ whose restriction
with the $i$-th face of $\bLambda_{n+2}^{n+2}$ agrees
with $y_i$, for every $i=1,\dots,n+2$. Since $A$ is
fibrant, we may then find $c\in A[n+2]$ such that
$\partial_ic=y_i$ for every $i=1,\dots,n+2$.
Set $b:=\partial_0c$; then
$$
\partial_ib=\partial_i\partial_0c=
\partial_0\partial_{i+1}c=\partial_0x_i
\qquad
\text{for every $i=0,\dots,n+1$}
$$
whence the assertion.

(b)$\Rightarrow$(a): Consider the map $f:[n]\times[1]\to[n+1]$
in the category $\Poset$ given by the rule :
$$
(i,0)\mapsto i
\quad\text{and}\quad
(i,1)\mapsto n+1
\quad
\text{for every $i\in[n]$}
$$
and let $\bar b:\bDelta_{n+1}\to A$ be the morphism
in $s.\Set$ corresponding to $b$; then it is easily
seen that $\bar b\circ\bDelta_f:\bDelta_n\times\bDelta_1\to A$
is a homotopy from $a$ to $\xi$ relative to $\partial\bDelta_n$
(details left to the reader).
\end{proof}

\sset\subsubsection{}\label{subsec_group-law-on-pi}
Let $(A,\xi)$ be any pointed simplicial set, such that
$A$ is fibrant. To justify the name of $\pi_n(A,\xi)$, we
shall endow it with a natural group structure, for every
$n>0$. To this aim, say that $\alpha,\beta\in\pi_n(A,\xi)$
are any two classes, and pick representatives
$x,y\in A(\xi,n)$ for $\alpha$ and respectively $\beta$.
Then, obviously the sequence
$$
\xi,\dots,\xi,x,y
$$
of $n+1$ elements of $A[n]$ satisfies the compatibility
condition of remark \ref{rem_fibrant-Kan}, with $k:=n$.
We may therefore find $z\in A[n+1]$ such that
$\partial_{n-1}z=x$, $\partial_{n+1}z=y$, and
$\partial_iz=\xi$ for $i=0,\dots,n-2$. The first
observation is :

\begin{lemma}\label{lem_group-law-on-pi}
With the notation of \eqref{subsec_group-law-on-pi},
the class of\/ $\partial_nz$ in $\pi_n(A,\xi)$ depends only
on $\alpha$ and $\beta$ (and not on $z$, nor the choice of
representatives $x$ and $y$).
\end{lemma}
\begin{proof} Let $x'$ and $y'$ be two other elements of
$A(\xi,n)$, and $h$ (resp. $h'$) a homotopy from $x$ to $x'$
(resp. from $y$ to $y'$) relative to $\partial\bDelta_n$.
Let also $z$ (resp. $z'$) be a morphism $\bDelta_{n+1}\to A$
such that
$$
z\circ\eps_i=\xi=z'\circ\eps_i
\qquad
\text{for $i=0,\dots,n-2$}
$$
and
$$
z\circ\eps_{n-1}=x
\qquad
z'\circ\eps_{n-1}=x'
\qquad
z\circ\eps_{n+1}=y
\qquad
z'\circ\eps_{n+1}=y'.
$$
Then there exists a unique morphism
$\psi:(\bDelta_{n+1}\times\partial\bDelta_1)\cup
(\bLambda^n_{n+1}\times\bDelta_1)\to A$ such that
\begin{itemize}
\item
$\psi\circ(\one_{\bDelta_{n+1}}\times\eps_0)=z$ and
$\psi\circ(\one_{\bDelta_{n+1}}\times\eps_1)=z'$.
\item
$\psi\circ(\eps_i\times\one_{\bDelta_1})$ factors
through $\xi$ for $i=0,\dots,n-2$.
\item
$\psi\circ(\eps_{n-1}\times\one_{\bDelta_1})=h$ and
$\psi\circ(\eps_{n+1}\times\one_{\bDelta_1})=h'$
\end{itemize}
and it is easily seen that $\psi\circ(\eps_n\times\one_{\bDelta_1})$
is a homotopy from $\partial_nz$ to $\partial_nz'$.
\end{proof}

\sset\subsubsection{}
Keep the notation of \eqref{subsec_group-law-on-pi}; lemma
\ref{lem_group-law-on-pi} says that the rule :
$(\alpha,\beta)\mapsto\partial_nz$ yields a well defined
pairing
\set\begin{equation}\label{eq_group-law-on-pi}
\pi_n(A,\xi)\times\pi_n(A,\xi)\to\pi_n(A,\xi)
\quad
(\alpha,\beta)\mapsto\alpha\cdot\beta
\qquad
\text{for every integer $n>0$}.
\end{equation}
For any $x\in A(\xi,n)$, we shall write $[x]$ for the class
of $x$ in $\pi_n(A,\xi)$.

\begin{theorem}
For every integer $n>0$, the pairing \eqref{eq_group-law-on-pi}
defines a group law on $\pi_n(A,\xi)$, whose neutral element is
the class $[\xi]$.
\end{theorem}
\begin{proof} Fix $n>0$ and any $\alpha\in\pi_n(A,\xi)$; to
see that $[\xi]\cdot\alpha=\alpha\cdot[\xi]=\alpha$,
pick any representative $x:\bDelta_n\to A$ for $\alpha$, and
set $z:=x\circ\eta_n$, $z':=x\circ\eta_{n+1}$. From the simplicial
identities \eqref{subsec_do-faces} it is easily seen that
$\partial_iz=\partial_iz'=\xi$ for $i=0,\dots,n-2$,
and $\partial_{n-1}z=\partial_{n+1}z'=x$, so $\partial_nz$
and $\partial_nz'$ represent respectively $[\xi]\cdot\alpha$
and $\alpha\cdot[\xi]$. But we have $\partial_nz=\partial_nz'=x$,
whence the claim.

Next, we check that the map
$$
\pi_n(A,\xi)\to\pi_n(A,\xi)
\quad:\quad
\beta\mapsto\alpha\cdot\beta
$$
is surjective. Indeed, given any $y\in A(\xi,n)$, there
exists a morphism $\phi:\bLambda_{n+1}^{n+1}\to A$ such that
$\phi\circ\eps_i=\xi$ for $i=0,\dots,n-2$, $\phi\circ\eps_{n-1}=x$
and $\phi\circ\eps_n=y$. We can then extend $\phi$ to a
morphism $z:\bDelta_n\to A$, and if we denote by $\beta$
(resp. $\gamma$) the class of $y$ (resp. of $\partial_{n+1}z$)
in $\pi_n(A,\xi)$, we see that $\alpha\cdot\gamma=\beta$.

Likewise, one shows that right multiplication by $\alpha$
defines a surjection on $\pi_n(A,\xi)$ (details left to the
reader). It remains only to check the associativity of
\eqref{eq_group-law-on-pi}, and to this aim, let
$x,y,z\in A(\xi,n)$ be any three elements, representing
respectively $\alpha,\beta,\gamma\in\pi_n(A,\xi)$.
Then there exist $u_{n-1},u_{n+1},u_{n+2}\in A[n+1]$
such that
$$
\partial_iu_{n-1}=\partial_iu_{n+1}=\partial_iu_{n+2}=\xi
\qquad
\text{for every $i=0,\dots,n-2$}
$$
and
$$
\begin{aligned}
\partial_{n-1}u_{n-1}=\,& x
& \qquad
\partial_{n+1}u_{n-1}=\,& y \\
\partial_{n-1}u_{n+1}=\,& \partial_nu_{n-1}
& \qquad
\partial_{n+1}u_{n+1}=\,& z \\
\partial_{n-1}u_{n+2}=\,& y
& \qquad
\partial_{n+1}u_{n+2}=\,& z.
\end{aligned}
$$
By remark \ref{rem_fibrant-Kan}, we may then find some
$v\in A[n+2]$ such that $\partial_iv$ equals $\xi$ for
$i=0,\dots,n-2$, and equals $u_i$ for $i=n-1,n+1,n+2$.
It follows that
$$
\partial_{n-1}\partial_nv=x
\qquad
\partial_{n+1}\partial_nv=\partial_nu_{n+2}
\quad\text{and}\quad
\partial_i\partial_nv=\xi
\quad
\text{for $i=0,\dots,n-2$}
$$
We may then compute :
$$
\begin{aligned}
(\alpha\cdot\beta)\cdot\gamma=\, &
[\partial_nu_{n-1}]\cdot\gamma \\
=\, & [\partial_nu_{n+1}] \\
=\, & [\partial_n\partial_{n+1}v] \\
=\, & [\partial_n\partial_nv] \\
=\, & \alpha\cdot(\beta\cdot\gamma)
\end{aligned}
$$
as required.
\end{proof}

\begin{remark}\label{rem_functorial-pi_n}
(i)\ \
Let us denote
$$
s.\Set^\mathrm{f}_\circ
$$
the category of {\em fibrant pointed simplicial sets}, whose
objects are all the pairs $(A,\xi)$ consisting of a simplicial
set $A$ and an element $\xi\in A[0]$. A morphism
$f:(A,\xi)\to(A',\xi')$ in $s.\Set^\mathrm{f}_\circ$ is
a morphism $f:A\to A'$ of simplicial sets such that
$f[0](\xi)=\xi'$. It follows easily from remark
\ref{rem_constant-homot}(ii) that, for every integer
$n>0$ (resp. for $n=0$), the rule
$(A,\xi)\mapsto\pi_n(A,\xi)$ (resp. $(A,\xi)\mapsto\pi_0(A)$)
yields a functor
$$
\pi_n:s.\Set^\mathrm{f}_\circ\to\Grp
\qquad
\text{(resp. $\pi_0:s.\Set^\mathrm{f}_\circ\to\Set_\circ$)}
$$
with values in the category of groups (resp. of pointed sets).
Namely, if $f$ is a morphism in $s.\Set^\mathrm{f}_\circ$ as in
the foregoing, and $x\in A(\xi,n)$ (resp. $x\in A[0]$) is any
element, we let $\pi_n(f,\xi)[x]$ (resp. $\pi_0(f)[x]$) be the
class of $f[n](x)$ in $\pi_n(A',\xi')$ (resp. in $\pi_0(A')$).
A simple inspection shows that this rule is well defined and
does yield a group homomorphism $\pi_n(A,\xi)\to\pi_n(A',\xi')$,
whenever $n>0$.

(ii)\ \
Let $X$ and $Y$ be two fibrant simplicial sets, $\xi\in X[0]$
any element, and denote also by $\xi$ the image of the
corresponding morphism $\bDelta_0\to X$, which is therefore a
simplicial subset of $X$. Let $f,g:X\to Y$ be two morphisms of
simplicial sets such that $f[0](\xi)=g[0](\xi)$ and
$f\sim_\xi g$. Then we claim that
$$
\pi_0(f)=\pi_0(g)
\quad\text{and}\quad
\pi_n(f,\xi)=\pi_n(g,\xi)
\qquad
\text{for every $n>0$}.
$$
Indeed, let $u$ be a homotopy from $f$ to $g$ relative to $\xi$;
if $x\in X(\xi,n)$, denote also by $x:\bDelta_n\to X$ the
corresponding morphism, and notice that
$u\circ(x\times\one_{\bDelta_1}):\bDelta_n\times\bDelta_1\to Y$
is a homotopy from $f[n](x)$ to $g[n](x)$ relative to
$\partial\bDelta_n$, whence the contention.
\end{remark}

\sset\subsubsection{}\label{subsec_telemarketing}
Let $p:X\to Y$ be a fibration between fibrant simplicial
sets $X$ and $Y$. Let also $\xi\in X[0]$ be any vertex,
and set $\xi':=p[0](\xi)$ and $F:=p^{-1}(\xi')$ (notation
of example \ref{ex_simplicial-sets}(vi)), so that $F$
is fibrant as well, by lemma \ref{lem_satur-RLP}(iv). For
any $n\in\N$ and any $a\in Y(\xi',n)$, we get a commutative
diagram
$$
\xymatrix{
\bLambda^0_n \ar[r]^-c \ar[d] & X \ar[d]^p \\
\bDelta_n \ar[r]^-a & Y
}$$
where $c$ is the unique morphism in $s.\Set$ that factors
through $\xi:\bDelta_0\to X$. We may then find $b\in X[n]$
such that $p[n](b)=a$ and $\partial_ib=\xi$ for every
$i=1,\dots,n$.

\begin{lemma}\label{lem_telemarketing}
In the situation of \eqref{subsec_telemarketing}, we
have $\partial_0b\in F(\xi,n-1)$, and the class
$[\partial_0b]\in\pi_{n-1}(F,\xi)$ depends only on
$[a]\in\pi_n(Y,\xi')$ (and is independent of the
representative $a$ for $[a]$ and of the simplex $b$).
\end{lemma}
\begin{proof} We have
$\partial_i\partial_0b=\partial_0\partial_{i+1}b=\xi$
for every $i=0,\dots,n-1$, whence the first assertion.
Next, suppose that $a'\in Y(\xi',n)$ is another
simplex such that $a\sim_{\partial\bDelta_n}a'$, and pick
any homotopy $h:\bDelta_n\times\bDelta_1\to Y$ from $a$
to $a'$ relative to $\partial\bDelta_n$. Let also
$b'\in X[n]$ be any simplex such that $p(b')=a'$ and
$\partial_ib'=\xi$ for $i=1,\dots,n$. We obtain a
commutative diagram
$$
\xymatrix{ (\bDelta_n\times\partial\bDelta_1)\cup
(\bLambda^0_n\times\bDelta_1) \ar[d] \ar[rr]^-{c'} & &
X \ar[d]^p \\
\bDelta_n\times\bDelta_1 \ar[rr]^-h & & Y
}$$
where $c'$ is the unique morphism whose restriction
to $\bDelta_n\times\bLambda^0_1$ (resp. to
$\bDelta_n\times\bLambda^1_1$) agrees with $b$ (resp.
with $b'$) and whose restriction to
$\bLambda^0_n\times\bDelta_1$ is the unique morphism
that factors through $\xi:\bDelta_0\to X$. Since $p$
is a fibration, we may then find a morphism
$h':\bDelta_n\times\bDelta_1\to X$ which extends $c'$,
and such that $p\circ h'=h$ (corollary \ref{cor_saturate}).
Set $h'':=h'\circ(\eps_0\times\one_{\bDelta_1}):
\bDelta_{n-1}\times\bDelta_1\to X$. It is easily seen
that $h''$ is a homotopy from $\partial_0b$ to $\partial_0b'$
relative $\partial\bDelta_{n-1}$, whence the lemma.
\end{proof}

\sset\subsubsection{}\label{subsec_homotopy-seq-fib}
Keep the situation of \eqref{subsec_telemarketing};
it follows from lemma \ref{lem_telemarketing} that the
rule $[a]\mapsto[\partial_0b]$ yields a well defined map
$$
\delta_n:\pi_n(Y,\xi')\to\pi_{n-1}(F,\xi)
\qquad
\text{for every $n>0$}
$$
and it is easily seen that $\delta_n[\xi]=[\xi']$,
{\em i.e.} $\delta_n$ is a map of pointed sets. Let
also $i:F\to X$ be the natural monomorphism; there
follows, for every $n\in\N$, a natural sequence of
pointed sets
$$
\pi_{n+1}(X,\xi)\xrightarrow{\ \pi_{n+1}(p,\xi)\ }
\pi_{n+1}(Y,\xi')\xrightarrow{\ \delta_{n+1}\ }\pi_n(F,\xi)
\xrightarrow{\ \pi_n(i,\xi)\ }\pi_n(X,\xi)
\xrightarrow{\ \pi_n(p,\xi)\ }\pi_n(Y,\xi')
$$
called the {\em homotopy sequence in degree $n$} attached
to the fibration $p$ (and the vertex $\xi$).
Let $\underline 1$ denote the monoid with one element
(the initial object in the category of monoids); noticed
that the category $\Set_\circ$ can also be regarded
naturally as the category of pointed (left)
$\underline 1$-modules (see \eqref{subsec_pointed-modules}),
and therefore the homotopy sequence in every degree is
naturally a sequence of pointed $\underline 1$-modules.

\begin{theorem}\label{th_homotopy-seq-fib}
With the notation of \eqref{subsec_homotopy-seq-fib},
the following holds :
\begin{enumerate}
\item
The homotopy sequence attached to $p$ is exact in
every degree $n\in\N$, {\em i.e.} we have
$$
\Img\,\pi_{n+1}(p,\xi)=\Ker\,\delta_{n+1}
\qquad
\Img\,\delta_{n+1}=\Ker\,\pi_n(i,\xi)
\qquad
\Img\,\pi_n(i,\xi)=\Ker\,\pi_n(p,\xi).
$$
\item
$\delta_n$ is a group homomorphism for every $n>1$.
\end{enumerate}
\end{theorem}
\begin{proof}(See remark \ref{rem_point-and-complete}(iii)
for the notion of kernel of a morphism of pointed
$\underline 1$-modules.)

(ii): Let $a_{n-1},a_n,a_{n+1}\in Y(\xi',n)$ be three
elements, $b_{n-1},b_n,b_{n+1}\in X[n]$ and $\phi\in Y[n+1]$
four simplices such that
\begin{itemize}
\item
$p(b_i)=a_i$ and $\partial_jb_i=\xi$ for every
$i=n-1,n,n+1$ and $j=1,\dots,n$
\item
$\partial_i\phi=\xi'$ for $i=0,\dots,n-2$ and
$\partial_i\phi=a_i$ for $i=n-1,n,n+1$.
\end{itemize}
We deduce a commutative diagram
$$
\xymatrix{ \bLambda_{n+1}^0 \ar[r]^-c \ar[d] & X \ar[d]^p \\
\bDelta_{n+1} \ar[r]^-\phi & Y
}$$
where $c$ is the unique morphism whose restriction to the
$i$-th face of $\bLambda^{n+1}_0$ factors through $\xi$ if
$i=1,\dots,n-2$, and equals $b_i$ for $i=n-1,n,n+1$.
Since $p$ is a fibration, we may then find
$\psi:\bDelta_{n+1}\to X$ whose restriction to
$\bLambda^0_{n+1}$ agrees with $c$, and such that
$p\circ\psi=\phi$. Set $t:=\partial_0\psi\in X[n]$; then
$p[n](t)=\partial_0(p[n+1](\psi))=\partial_0\phi=\xi$,
so that $t\in F[n]$. Moreover,
$\partial_it=\partial_0\partial_{i+1}\psi$ for every
$i=0,\dots,n$, hence $\partial_it=\xi$ for $i<n-2$
and $\partial_it=\partial_0b_{i+1}$ for $i=n-2,n-1,n$.
This means that
$$
[\partial_0b_n]=[\partial_0b_{n-1}]\cdot[\partial_0b_{n+1}]
\qquad
\text{in $\pi_{n-1}(F,\xi)$}
$$
and on the other hand we have $[a_n]=[a_{n-1}]\cdot[a_{n+1}]$
in $\pi_n(Y,\xi')$, whence the assertion.

(i): The inclusion $\Img\,\pi_n(i,\xi)\subset\Ker\,\pi_n(p,\xi)$
is immediate, and the identity
$\Img\,\delta_{n+1}=\Ker\,\pi_n(i,\xi)$ follows easily from
lemma \ref{lem_hom-groups}, by inspection of the definition
of $\delta_{n+1}$. The inclusion
$\Img\,\pi_{n+1}(p,\xi)\subset\Ker\,\delta_{n+1}$ is likewise
immediate, from the definition of $\delta_{n+1}$ (details left
to the reader).

$\bullet$\ \ 
Next, let us show that
$\Ker\,\delta_{n+1}\subset\Img\,\pi_{n+1}(p,\xi)$. Indeed,
let $a\in Y(\xi',n+1)$, pick $b\in X[n+1]$ such that
$p[n+1](b)=a$ and $\partial_ib=\xi$ for every $i=1,\dots,n+1$,
and suppose that $\partial_0b\sim_{\partial\bDelta_n}\xi$.
Let $h:\bDelta_n\times\bDelta_1\to F$ be a homotopy
from $\partial_0b$ to $\xi$ relative to $\partial\bDelta_n$.
Then we get a unique morphism
$$
t:(\bDelta_{n+1}\times\bLambda^1_1)\cup
(\partial\bDelta_{n+1}\times\bDelta_1)
\to X
$$
whose restriction to $\bDelta_{n+1}\times\bLambda^1_1$ agrees
with $b$, whose restriction to $\bLambda^0_{n+1}\times\bDelta_1$
factors through $\xi$, and whose restriction to
$(\Img\,\eps_0)\times\bDelta_1\subset
\partial\bDelta_{n+1}\times\bDelta_1$ agrees with $h$.
Since $X$ is fibrant, $t$ extends to a morphism
$t':\bDelta_{n+1}\times\bDelta_1\to X$, and we denote
by $x\in X[n+1]$ the restriction of $t'$ to
$\bDelta_{n+1}\times\bLambda^0_1$. Then it is easily
seen that $x\in X(\xi,n+1)$, and $p\circ t'$ is a homotopy
relative to $\partial\bDelta_{n+1}$ from $a$ to $p(x)$,
whence the assertion.

$\bullet$\ \ 
Lastly, let $b\in X(\xi,n)$ be any simplex such that
$[p(b)]=[\xi']$ in $\pi_n(Y,\xi')$, and pick any
homotopy $h:\bDelta_n\times\bDelta_1\to Y$ from $p(b)$
to $\xi'$ relative to $\partial\bDelta_n$. We get a
commutative diagram
$$
\xymatrix{
(\bDelta_n\times\bLambda^0_1)\cup
(\partial\bDelta_n\times\bDelta_1) \ar[rr]^-\phi \ar[d] & &
X \ar[d]^p \\
\bDelta_n\times\bDelta_1 \ar[rr]^-h & & Y
}$$
where $\phi$ is the unique morphism whose restriction
to $\partial\bDelta_n\times\bDelta_1$ factors through
$\xi$, and whose restriction to $\bDelta_n\times\bLambda^0_1$
agrees with $b$. Since $p$ is a fibration, we may then
find a morphism $h':\bDelta_n\times\bDelta_1\to X$
extending $\phi$ and such that $p\circ h'=h$. Set
$b':=h'\circ(\one_{\bDelta_n}\times\eps_1)$, so that
$h'$ is a homotopy from $b$ to $b'$ relative to
$\partial\bDelta_n$; however, notice that $p(b')=\xi'$,
so $b'\in F(\xi,n)$, therefore $[b]\in\pi_n(i,\xi)$,
which concludes the proof of the theorem.
\end{proof}

\sset\subsubsection{}\label{subsec_device-groupoid}
The following device shall be useful for explaining the
functorial behaviour of the homotopy groups under change
of base points. Let $A$ be any fibrant simplicial set; we set
$$
\pi(A):=A[1]/\!\sim_{\partial\bDelta_1}.
$$
Given $a\in A[1]$, we denote by $[a]\in\pi(A)$ the class
of $a$. Evidently, the vertices $\partial_ia\in A[0]$
($i=0,1$) depend only on $[a]$, so we shall also
denote them by $\partial_i[a]$. We use the notation
$$
[a]:x\to y
$$
to signal that $x=\partial_0[a]$ and $y=\partial_1[a]$.
If $[b]:y\to z$ is another class, we define a composition
$[b]\circ[a]:x\to z$ as follows. Let $\phi:\bLambda_2^1\to A$
be the unique morphism such that $\phi\circ\eps_0=a$ and
$\phi\circ\eps_2=b$; since $A$ is fibrant, $\phi$ extends to
a morphism $c:\bDelta_2\to A$, and we let
$$
[b]\circ[a]:=[\partial_1c].
$$
In order to show that this rule is well defined, suppose
that $a',b'\in A[1]$ and $c'\in A[2]$ are any three other
simplices such that $[a]=[a']$, $[b]=[b']$, $\partial_0c'=a'$
and $\partial_2c'=b'$. Pick a homotopy $h$ (resp. $h'$)
from $a$ to $a'$ (resp. from $b$ to $b'$) relative to
$\partial\bDelta_1$; there follows a morphism
$$
t:(\bDelta_2\times\partial\bDelta_1)\cup
(\bLambda_2^1\times\bDelta_1)\to A
$$
whose restriction to $\bDelta_2\times\bLambda^0_1$
(resp. to $\bDelta_2\times\bLambda^1_1$) agrees with
$c$ (resp. with $c'$) and whose restriction to
$(\Img\,\eps_0)\times\bDelta_1$ (resp. to
$(\Img\,\eps_2)\times\bDelta_1$) agrees with $h$ (resp.
with $h'$). Since $A$ is fibrant, $t$ extends to a
morphism $\bDelta_2\times\bDelta_1\to A$, whose
restriction to $(\Img\,\eps_1)\times\bDelta_1$ is a
homotopy from $\partial_1c$ to $\partial_1c'$ relative
to $\partial\bDelta_1$, whence the claim. Moreover,
we claim that this composition law is associative :
indeed, the proof is the same as that of the associativity
of the group law on $\pi_1(A,\xi)$ (details left to
the reader). Furthermore, for every vertex $x\in A[0]$,
let $\one_x:\bDelta_1\to A$ be the unique morphism which
factors through $x$; then $[a]\circ[\one_x]=[a]$ for
every $y\in A[0]$ and every $[a]:x\to y$, since
$$
\partial_0(a\circ\eta_1)=\one_x
\qquad
\partial_1(a\circ\eta_1)=a
\qquad
\partial_2(a\circ\eta_1)=a.
$$
Likewise, we see that $[\one_x]\circ[b]=[b]$ for every
$[b]:y\to x$, by considering $\partial(b\circ\eta_0)=b$
(details left to the reader). Lastly, for any $[a]:x\to y$
consider the unique morphism $t:\bLambda^0_2\to A$ such that
$t\circ\eps_1=\one_x$ and $t\circ\eps_0=a$; since $A$
is fibrant, $t$ extends to a morphism $c:\bDelta_2\to A$,
and it follows that $[\partial_2c]\circ[a]=\one_x$.
Similarly, we see easily that $[a]$ admits a right
inverse. We conclude that the datum
$$
\underline\pi(A):=(A[0],\pi(A),\circ)
$$
is a {\em groupoid}, {\em i.e.} a category all whose
morphisms are isomorphisms. It shall be called the
{\em fundamental groupoid of $A$}.

\begin{remark}\label{rem_fundam-groupoid}
Let $f:X\to Y$ be any morphism of fibrant simplicial sets.
A simple inspection of the definitions shows that $f$
induces a functor
$$
\underline\pi(f):\underline\pi(X)\to\underline\pi(Y).
$$
Explicitly, $\underline\pi(f)$ is the functor whose map
on objects $X[0]\to Y[0]$ is given by $f[0]$, and whose
map on morphisms is induced by $f[1]$. It is then clear
that the rule $f\mapsto\underline\pi(f)$ defines a
functor
$$
\underline\pi:s.\Set\to\Gpd
$$
with values in the full subcategory of $\bCat$ whose
objects are all (small) groupoids.
\end{remark}

\sset\subsubsection{}\label{subsec_groupoid-action}
Now, let $n\in\N$ be any integer, $x\in A[0]$ any
vertex, $b:\bDelta_n\to A$ any element of $A(x,n)$,
and $[\alpha]:x\to y$ any element of $\pi(A)$.
We have a unique morphism
$$
t:(\partial\bDelta_n\times\bDelta_1)\cup
(\bDelta_n\times\bLambda^1_1)\to A
$$
such that :
\begin{itemize}
\item
the restriction of $t$ to $\partial\bDelta_n\times\bDelta_1$
equals $\alpha\circ p$, where
$p:\partial\bDelta_n\times\bDelta_1\to\bDelta_1$
is the natural projection
\item
the restriction of $t$ to $\bDelta_n\times\bLambda^1_1$
agrees with $b$.
\end{itemize}
Since $A$ is fibrant, $t$ extends to a morphism
$h:\bDelta_n\times\bDelta_1\to A$, and we denote by
$\partial_1h\in A[n]$ the restriction of $h$ to
$\bDelta_n\times\bLambda^0_1$. Notice that
$\partial_1h\in A(y,n)$.

\begin{lemma}\label{lem_fundamental-action}
With the notation of \eqref{subsec_groupoid-action}, the
class $[\partial_1h]\in\pi_n(A,y)$ depends only on $[\alpha]$
and $[b]\in\pi_n(A,x)$ (and is independent of $h$ and of the
representatives $\alpha$ and $b$ for these classes).
\end{lemma}
\begin{proof} Say that $[\beta]=[\alpha]$ in $\pi(A)$, so
we have a homotopy $\gamma:\bDelta_1\times\bDelta_1\to A$
from $\alpha$ to $\beta$ relative to $\partial\bDelta_1$.
Suppose also that $[b]=[b']$ in $\pi_n(A,x)$, and let
$H:\bDelta_n\times\bDelta_1\to A$ be a given homotopy
from $b$ to $b'$ relative to $\partial\bDelta_n$.
Pick also a morphism $h':\bDelta_n\times\bDelta_1\to A$
whose restriction to $\partial\bDelta_n\times\bDelta_1$
(resp. to $\bDelta_n\times\bLambda^1_1$) equals
$\beta\circ p$ (resp. $b'$). To ease notation, set
$$
K:=(\bLambda^1_1\times\bDelta_1)\cup
(\bDelta_1\times\partial\bDelta_1)
\qquad\text{and}\qquad
L:=\bDelta_1\times\bDelta_1
$$
and notice that the monomorphism $K\to L$ is anodyne.
We then have a unique morphism
$$
u:(\partial\bDelta_n\times L)\cup(\bDelta_n\times K)\to A
$$
such that
\begin{itemize}
\item
the restriction of $u$ to $\partial\bDelta_n\times L$ equals
$q\circ\gamma$, where $q:\partial\bDelta_n\times L\to L$ is
the natural projection
\item
the restriction of $u$ to
$\bDelta_n\times(\bLambda^1_1\times\bDelta_1)$ equals $H$
\item
the restriction of $u$ to
$\bDelta_n\times(\bLambda^1_1\times\bDelta_1)$ (resp. to
$\bDelta_n\times(\bLambda^0_1\times\bDelta_1)$) equals
$h$ (resp. $h'$).
\end{itemize}
By corollary \ref{cor_saturate}, the morphism $u$
extends to a morphism $h'':\bDelta_n\times L\to A$,
and it is easily seen that the restriction of $h''$
to $\bDelta_n\times(\bLambda^0_1\times\bDelta_1)$ is
a homotopy from $\partial_1h$ to $\partial_1h'$ relative
to $\partial\bDelta_n$, whence the assertion.
\end{proof}

\sset\subsubsection{}\label{subsec_fundamental-action}
In light of lemma \ref{lem_fundamental-action}, the
rule $(\alpha,b)\mapsto[\partial_1h]$ yields a well
defined mapping
$$
\pi_n(A,[\alpha]):\pi_n(A,x)\to\pi_n(A,y)
$$
for every $n\in\N$, every $x,y\in A[0]$ and every
$[\alpha]:x\to y$. For $n=0$, a simple inspection shows
that the resulting map
$\pi_0(A,[\alpha]):\pi_0(A)\to\pi_0(A)$ is the identity.
Thus, in the following we assume that $n>0$, in which
case we have :

\begin{proposition} With the notation of
\eqref{subsec_fundamental-action}, the following holds
for every $n>0$ :
\begin{enumerate}
\item
$\pi_n(A,[\alpha])$ is a group homomorphism, and
$\pi_n(A,[\one_x])=\one_{\pi_n(A,x)}$ for every $x\in A[0]$.
\item
If $[\beta]:y\to z$ is any other element of $\pi(A)$,
we have
$$
\pi_n(A,[\beta])\circ\pi_n(A,[\alpha])=
\pi_n(A,[\beta]\circ[\alpha]).
$$
\item
Therefore, the rule
$$
x\mapsto\pi_n(A,x)
\qquad
[\alpha]\mapsto\pi_n(A,[\alpha])
\qquad
\text{for every $x\in A[0]$ and every $[\alpha]\in\pi(A)$}
$$
yields a well defined functor
$$
\pi_n:\underline\pi(A)\to\Grp.
$$
\end{enumerate}
\end{proposition}
\begin{proof}(i): Let $b,b'\in A(x,n)$ be any two
elements, and pick a morphism $c:\bDelta_{n+1}\to A$
such that $\partial_ic=x$ for $i=0,\dots,n-2$,
$\partial_{n-1}c=b$ and $\partial_{n+1}c=b'$. Let
also $[\alpha]:x\to y$ be any element of $\pi(A)$, and
choose morphisms $h,h':\bDelta_n\times\bDelta_1\to A$
whose restrictions to $\partial\bDelta_n\times\bDelta_1$
equal $p\circ\alpha$ (where $p$ is as in
\eqref{subsec_groupoid-action}) and whose restrictions
to $\bDelta_n\times\bLambda^1_1$ agree with $b$ and
respectively $b'$. Then there exists a unique morphism
$T:(\partial\bDelta_{n+1}\times\bDelta_1)\cup
(\bDelta_{n+1}\times\bLambda^1_1)\to A$ such that
\begin{itemize}
\item
for $i=0,\dots,n-2$, the restriction of $T$ to
$(\Img\,\eps_i)\times\bDelta_1$ equals $p'\circ\alpha$,
where $p':(\Img\,\eps_i)\times\bDelta_1\to\bDelta_1$ is
the projection
\item
the restriction of $T$ to $(\Img\,\eps_{n-1})\times\bDelta_1$
(resp. to $(\Img\,\eps_{n+1})\times\bDelta_1$) agrees with
$h$ (resp. with $h'$)
\item
the restriction of $T$ to $\bDelta_{n+1}\times\bLambda^1_1$
agrees with $c$.
\end{itemize}
Since $A$ is fibrant, $T$ extends to a morphism
$H:\bDelta_{n+1}\times\bDelta_1\to A$, and we denote by
$U:\bDelta_n\times\bDelta_1\to A$ the restriction of $H$
to $(\Img\,\eps_n)\times\bDelta_1$. By inspecting the
definition, it is easily seen that $\partial_1U$ is a
representative for both $\pi_n(A,[\alpha])([b]\cdot[b'])$
and $\pi_n(A,[\alpha])([b])\cdot\pi_n(A,[\alpha])([b'])$,
so $\pi_n(A,[\alpha])$ is a group homomorphism. A simple
inspection shows that $\pi_n(A,[\one_x])$ is the identity
map of $\pi_n(A,x)$.

(ii): Pick any morphism $\phi:\bDelta_2\to A$
such that $\partial_0\phi=\alpha$ and $\partial_2\phi=\beta$.
Let also $b$ and $h$ as in \eqref{subsec_groupoid-action}, and
choose similarly a morphism $h':\bDelta_n\times\bDelta_1\to A$
whose restriction to $\bDelta_n\times\bLambda^1_1$ equals
$\partial_1h$, and whose restriction to
$\partial\bDelta_n\times\bDelta_1$ equals $p\circ\beta$.
Then there exists a unique morphism
$u:(\partial\bDelta_n\times\bDelta_2)\cup
(\bDelta_n\times\bLambda^1_2)\to A$ such that
\begin{itemize}
\item
the restriction of $u$ to $\partial\bDelta_n\times\bDelta_2$
equals $q\circ\phi$, where
$q:\partial\bDelta_n\times\bDelta_2\to\bDelta_2$ is the
projection
\item
the restriction of $u$ to $\bDelta_n\times(\Img\,\eps_0)$
(resp. to $\bDelta_n\times(\Img\,\eps_2)$) equals $h$
(resp. $h'$).
\end{itemize}
Since $A$ is fibrant, $u$ extends to a morphism
$w:\bDelta_n\times\bDelta_2\to A$, and we denote by
$w_1:\bDelta_n\times\bDelta_1\to A$ the
restriction of $w$ to $\bDelta_n\times(\Img\,\eps_1)$.
Hence, the restriction of $w_1$ to
$\partial\bDelta_n\times\bDelta_1$ equals
$p\circ(\partial_1\phi)$, and its restriction
$\partial_1w_1$ to $\bDelta_n\times\bLambda^0_1$
(resp. to $\bDelta_n\times\bLambda^1_1$) equals $b$
(resp. $\partial_1h'$). Unwinding the definitions,
we find that
$$
\begin{aligned}
\pi_n(A,[\beta])\circ\pi_n(A,[\alpha])([b])= &\,
\pi_n(A,[\beta])([\partial_1h]) \\
= &\, [\partial_1w_1] \\
= &\, \pi_n(A,[\partial_1\phi])([b]) \\
= &\, \pi_n(A,[\beta]\circ[\alpha])([b])
\end{aligned}
$$
whence the contention. Of course, (iii) results by
combining (i) and (ii).
\end{proof}

\begin{remark}\label{rem_sothat}
Let $f:X\to Y$ be any morphism of fibrant simplicial sets,
$[\alpha]:x\to x'$ a morphism in $\underline\pi(X)$, and
set $y:=f[0](x)$, $y':=f[0](y')$ and $\beta:=f[1](\alpha)$,
so that
$[\beta]=\underline\pi(f)[\alpha]$ in $\underline\pi(Y)$
(notation of remark \ref{rem_fundam-groupoid}). A direct
inspection shows that the resulting diagram
$$
\xymatrix{
\pi_n(X,x) \ar[rr]^-{\pi_n(X,[\alpha])} \ar[d]_{\pi_n(f,x)} & &
\pi_n(X,x') \ar[d]^{\pi_n(f,x')} \\
\pi_n(Y,y) \ar[rr]^-{\pi_n(Y,[\beta])} & & \pi_n(Y,y')
}$$
commutes (details left to the reader).
\end{remark}

\begin{definition}
Consider any commutative diagram of simplicial sets :
$$
\xymatrix{ X' \ar@<.5ex>[rr]^-f \ar@<-.5ex>[rr]_-g
\ar[rd]_{p'} & & X \ar[ld]^p \\
& Y.
}$$
\begin{enumerate}
\item
A {\em $p$-fibrewise homotopy from $f$ to $g$} is
a homotopy $h$ from $f$ to $g$ such that $p\circ h$
is a constant homotopy (see remark \ref{rem_constant-homot}(i)).
\item
Let also $K\subset X'$ be any simplicial subset. A
{\em $p$-fibrewise homotopy from $f$ to $g$ relative to $K$}
is a $p$-fibrewise homotopy from $f$ to $g$ which is also a
homotopy relative to $K$ (see definition \ref{def_homotop_rel}).
We write
$$
f\overset{p}{\sim} g
\qquad
\text{(resp. $f\overset{p}{\sim}_Kg$)}
$$
if there exists a $p$-fibrewise homotopy from $f$ to
$g$ (resp. relative to $K$).
\item
Recall that, for every $n\in\N$, we may regard any
$x\in X[n]$ as a morphism $x:\bDelta_n\to X$. Then,
we say that $p$ is a {\em minimal fibration}, if $p$
is a fibration, and
$$
a=b
\quad\Leftrightarrow\quad
a\overset{p}{\sim}_{\partial\bDelta_n}b
$$
for every $n\in\N$ and every $a,b\in X[n]$ such that
$p[n](a)=p[n](b)$.
\item
A {\em strong deformation retract} of $p$ is a datum
$(i,r,h)$ consisting of a monomorphism $i:X'\to X$,
an epimorphism $r:X\to X'$ such that
$$
p\circ i\circ r=p
\qquad\text{and}\qquad
r\circ i=\one_{X'}
$$
and a $p$-fibrewise homotopy $h:X\times\bDelta_1\to X$
from $\one_X$ to $i\circ r$ relative to $X'$.
\item
We say that $f$ is a {\em fibrewise homotopy equivalence
from $p'$ to $p$}, if there exists a morphism of simplicial
sets $f':X\to X'$ with
$$
p'\circ f'=p
\qquad
f'\circ f\overset{p'}{\sim}\one_{X'}
\qquad
f\circ f'\overset{p}{\sim}\one_X.
$$
\item
For any simplicial set $A$, let $p_A:A\to\bDelta_0$ be the
unique morphism of simplicial sets. If $A$ is fibrant, we
say that $A$ is {\em minimal}, if $p_A$ is a minimal fibration.
We say that a morphism of simplicial sets $f:A\to A'$ is
a {\em homotopy equivalence}, if $f$ is a fibrewise
homotopy equivalence from $p_A$ to $p_{A'}$. We say that
a monomorphism $i:A'\to A$ is a
{\em strong deformation retract of $A$}, if there exist
a morphism $r:A\to A'$ and a homotopy $h$ from $\one_A$
to $i\circ r$, such that $(i,r,h)$ is a strong deformation
retract of $p_A$.
\end{enumerate}
\end{definition}

\begin{remark}\label{rem_equivalences}
(i)\ \
By inspecting the definitions, it is easily seen that
a composition of fibrewise homotopy equivalences is a
fibrewise homotopy equivalence. Especially, a composition
of homotopy equivalences is a homotopy equivalence.

(ii)\ \
Let $(i,r,h)$ be any strong deformation retract; then,
clearly, both $i$ and $r$ are fibrewise homotopy equivalences.
\end{remark}

\begin{lemma}\label{lem_fibrewise-equiv}
Let $p:X\to Y$ be any fibration, $p':X'\to Y$ any object
of $s.\Set/Y$, and $i:K\to X'$ any monomorphism.
Then :
\begin{enumerate}
\item
$\overset{p}{\sim}_K$ is an equivalence relation
on $\Hom_{s.\Set/Y}(p',p)$.
\item
More precisely, let $f,f',f''\in\Hom_{s.\Set/Y}(p',p)$
be any three elements such that
$$
f\overset{p}{\sim}f'
\qquad
f'\overset{p}{\sim}_Kf''
$$
and $h$ a $p$-fibrewise homotopy from $f$ to $f'$.
Then there exists a $p$-fibrewise homotopy $h'$ from
$f$ to $f''$ such that
$$
h'\circ(i\times\one_{\bDelta_1})=h\circ(i\times\one_{\bDelta_1}).
$$
\item
Let $f,f'$ be as in {\em(ii)}, and suppose that
$f$ is a fibrewise homotopy equivalence. Then the
same holds for $f'$.
\end{enumerate}
\end{lemma}
\begin{proof}(i): It is a relative variant of the proof
of theorem \ref{th_hom-is-equiv}. Namely, consider first
the case of $\overset{p}{\sim}_K$; on the one hand, by
proposition \ref{prop_fibrant-hom}, the morphism
$$
t:s.\cHom(X',X)\to
P:=s.\cHom(X',Y)\times_{s.\cHom(K,Y)}s.\cHom(K,X)
$$
induced by $p$ and $i$ is a fibration. On the other hand,
let $f,f'\in\Hom_{s.\Set/Y}(p',p)$ be any two elements; the
datum of a homotopy $h$ from $f$ to $f'$ relative to $K$
is the same as that of a morphism
$$
h^*:\bDelta_1\to s.\cHom(X',X)
$$
such that
$s.\cHom(i,X)\circ h^*$ factors through a morphism
$a:\bDelta_0\to s.\cHom(K,X)$ and the unique morphism
$\pi:\bDelta_1\to\bDelta_0$. Then, $h$ is a $p$-fibrewise
homotopy if and only if
$$
\begin{aligned}
s.\cHom(K,p)\circ a=\, & s.\cHom(i,Y)\circ b \\
s.\cHom(X',p)\circ h^*=\, & b\circ\pi
\end{aligned}
$$
for the unique morphism $b:\bDelta_0\to s.\cHom(X',Y)$
corresponding to $p'$.
Now, the pair $(a,b)$ determines a unique morphism
$(a,b):\bDelta_0\to P$, and the foregoing conditions
tell us that $t\circ h^*=(a,b)\circ\pi$. Thus, set
$Q:=t^{-1}(a,b)$ (notation of example
\ref{ex_simplicial-sets}(vi)); by lemma
\ref{lem_satur-RLP}(iv), the simplicial set $Q$ is
fibrant, and we have a natural bijection
$$
Q[0]\isom\{f\in\Hom_{s.\Set/Y}(p',p)~|~f\circ i=a\}
$$
which identifies the restriction of the relation
$\overset{p}{\sim}_K$ with the relation $\sim$ on
$Q[0]$. The latter is an equivalence relation, by
lemma \ref{lem_homotopy-vertex}, whence the contention.

(ii): Let again $b:\bDelta_0\to s.\cHom(X',Y)$ be
the morphism corresponding to $p'$, and set
$$
P_b:=q^{-1}(b)
\qquad
Q_b:=P_b\times_P s.\cHom(X',X)
$$
where $q:P\to s.\cHom(X',Y)$ is the natural projection.
So, the morphism $t$ restricts to a fibration $t_b:Q_b\to P_b$.
By assumption, $f'$ and $f''$ can be regarded as morphisms
$\bDelta_0\to s.\cHom(X',X)$ such that
$a:=s.\cHom(i,X)\circ f'=s.\cHom(i,X)\circ f''$,
and there exists a morphism
$g:\bDelta_1\to s.\cHom(X',X)$ such that
$s.\cHom(i,X)\circ g$ factors through $a$ and
$s.\cHom(X',p)\circ g$ factors through $b$,
and moreover $g\circ\eps_1=f'$, $g\circ\eps_0=f''$.
Likewise, the homotopy $h$ can be regarded as a morphism
$\bDelta_1\to s.\cHom(X',X)$ such that
$s.\cHom(X',p)\circ h$ factors through $b$,
and moreover $h\circ\eps_1=f$, $h\circ\eps_0=f'$.
The pair $(h,g)$ then amounts to a morphism
$(h,g):\bLambda_2^1\to Q_b$ fitting into a commutative
diagram
$$
\xymatrix{
\bLambda_2^1 \ar[d]_{\iota_2^1} \ar[rr]^-{(h,g)} & &
Q_b \ar[d]^{t_b} \\
\bDelta_2 \ar[r]^-{\eta_0} & \bDelta_1 \ar[r] & P_b
}$$
(namely, $(h,g)$ is the unique morphism such that
$(h,g)\circ\eps_0=h$ and $(h,g)\circ\eps_2=g$).
Thus, $(h,g)$ extends to a morphism
$u:\bDelta_2\to Q_b$ such that $t_b\circ u$ factors
as well through $\eta_0$; it is easily seen that
the morphism $u\circ\eps_1$ yields the sought
homotopy $h'$.

(iii): By assumption, there exists a morphism
$g:X\to X'$ of $s.\Set/Y$ such that
$g\circ f\overset{p}{\sim}\one_{X'}$ and
$f\circ g\overset{p'}{\sim}\one_X$. On the other
hand, since $f\overset{p}{\sim}f'$, remark
\ref{rem_constant-homot}(ii) implies that
$g\circ f\overset{p'}{\sim}g\circ f'$ and
$f\circ g\overset{p}{\sim}f'\circ g$. Then
the assertion follows from (i).
\end{proof}

\sset\subsubsection{}\label{subsec_pbk-is-strongly-h}
Consider a fibration $p:X\to Y$ of simplicial sets, and
two morphisms $f_0,f_1:A\to Y$, and for $i=0,1$ define
the simplicial set $X_i$ as the fibre product in the
cartesian diagram
$$
\xymatrix{ X_i \ar[r] \ar[d]_{p_i} & X \ar[d]^p \\
A \ar[r]^-{f_i} & Y.
}$$

\begin{proposition}
In the situation of \eqref{subsec_pbk-is-strongly-h},
suppose that $f_0\sim f_1$ (notation of
\eqref{subsec_homotop-rel}). Then there exists a
fibrewise homotopy equivalence from $p_0$ to $p_1$.
\end{proposition}
\begin{proof} Let $h:A\times\bDelta_1\to Y$ be a homotopy
from $f_0$ to $f_1$, and define the simplicial set $H$ as
the fibre product in the resulting diagram whose two
square subdiagrams are cartesian
$$
{\diagram X_i \ar[r]^{j_i} \ar[d]_{p_i} &
H \ar[r] \ar[d]_{p_H} & X \ar[d]^p \\
A \ar[r]^-{A\times\eps_i} & A\times\bDelta_1 \ar[r]^-h & Y
\enddiagram}
\qquad
\text{for $i=0,1$.}
$$
Then, $p_H$ is a fibration (lemma \ref{lem_satur-RLP}(iv)),
hence there exist commutative diagrams
$$
{\diagram
X_i \ar[rr]^-{j_i} \ar[d]_{X_i\times\eps_i} & &
H \ar[d]^{p_H} \\
X_i\times\bDelta_1 \ar[rr]^-{p_i\times\bDelta_1}
\ar[rru]^-{\vartheta_i} & &
A\times\bDelta_1
\enddiagram}
\qquad
\text{for $i=0,1$}
$$
(proposition \ref{prop_saturate}) whence two commutative
diagrams with cartesian square subdiagrams
$$
{\diagram X_i \ar[r]^-{\omega_i} \ar[d]_{X_i\times\eps_{1-i}}
\ar@/^2pc/[rr]^-{p_i} &
X_{1-i} \ar[r]^-{p_{1-i}} \ar[d]^{j_{1-i}} &
A \ar[d]^{A\times\eps_{1-i}} \\
X_i\times\bDelta_1 \ar[r]^-{\vartheta_i}
\ar@/_2pc/[rr]_-{p_i\times\bDelta_1} &
H \ar[r]^-{p_H} & A\times\bDelta_1
\enddiagram}
\qquad
\text{for $i=0,1$}.
$$
Notice now that
$$
\vartheta_0\circ(X_0\times\eps_1)=j_1\circ\omega_0=
\vartheta_1\circ(X_1\times\eps_1)\circ\omega_0=
\vartheta_1\circ(\omega_0\times\bDelta_1)\circ(X_0\times\eps_1)
$$
from which we deduce that there exists a commutative
diagram
$$
\xymatrix{
X_0\times\bLambda_2^0 \ar[rr]^-\beta \ar[d]_{X_0\times\iota_2^0}
& & H \ar[d]^{p_H} \\
X_0\times\bDelta_2 \ar[rr]^-{p_0\times\eta_1} \ar[rru]^-\gamma
& & A\times\bDelta_1
}$$
where $\beta$ is the unique morphism such that
$$
\beta\circ(X_0\times\eps_1)=\vartheta_0
\qquad\text{and}\qquad
\beta\circ(X_0\times\eps_2)=
\vartheta_1\circ(\omega_0\times\bDelta_1).
$$
Set $\gamma_0:=\gamma\circ(X_0\times\eps_0):
X_0\times\bDelta_1\to H$, and notice that
$\eta_1\circ\eps_0:\bDelta_1\to\bDelta_1$ factors through
$\eps_0:\bDelta_0\to\bDelta_1$; it follows that $\gamma_0$
factors through $j_0:X_0\to H$ and a morphism
$\bar\gamma_0:X_0\times\bDelta_1\to X_0$. 
Lastly, notice that
$$
\begin{aligned}
j_0\circ\bar\gamma_0\circ(X_0\times\eps_0)=\,&
\vartheta_0\circ(X_0\times\eps_0)=j_0 \\
j_0\circ\bar\gamma_0\circ(X_0\times\eps_1)=\,&
\vartheta_1\circ(\omega_0\times\bDelta_1)\circ(X_0\times\bDelta_0)
=\vartheta_1\circ(X_1\times\eps_0)\circ\omega_0=
j_0\circ\omega_1\circ\omega_0
\end{aligned}
$$
whence
$$
\bar\gamma_0\circ(X_0\times\eps_0)=\one_{X_0}
\qquad\text{and}\qquad
\bar\gamma_0\circ(X_0\times\eps_1)=\omega_1\circ\omega_0
$$
so $\bar\gamma_0$ is a homotopy from $\one_{X_0}$ to
$\omega_1\circ\omega_0$. A similar construction yields
likewise a homotopy from $\one_{X_1}$ to
$\omega_0\circ\omega_1$, whence the proposition (details
left to the reader).
\end{proof}

\begin{proposition}\label{prop_minimals-exist}
Every fibration $p:X\to Y$ admits a strong deformation
retract $(j,r,h)$ such that $p\circ j$ is a minimal fibration.
\end{proposition}
\begin{proof} For every $k,n\in\N$ with $n\geq k$,
denote by
$$
\Sk_kX\xrightarrow{\ \eps^{(k,n)}_X\ }
\Sk_nX\xrightarrow{\ \eta^{(n)}_X\ }X
$$
the inclusion maps (notation of remark \ref{rem_degenerate}(v))
and set $p^{(n)}:=p\circ\eta_X^{(n)}$ for every $n\in\N$.

\begin{claim}\label{cl_approx-def-retract}
For every $n\in\N$ there exist morphisms
$$
j^{(n)}:Z^{(n)}\to\Sk_nX
\qquad
r^{(n)}:\Sk_nX\to Z^{(n)}
\qquad
h^{(n)}:\Sk_nX\times\bDelta_1\to X
$$
such that :
\begin{enumerate}
\alphaenu
\item
$r^{(n)}\circ j^{(n)}=\one_{Z^{(n)}}$ and
$p^{(n)}\circ j^{(n)}\circ r^{(n)}=p^{(n)}$.
\item
$h^{(n)}$ is a homotopy from $\eta_X^{(n)}$ to
$\eta_X^{(n)}\circ j^{(n)}\circ r^{(n)}$.
\item
Both $p\circ h^{(n)}$ and the restriction
$ h^{(n)}_Z:Z^{(n)}\times\bDelta_1\to X$ of $h^{(n)}$ are
constant homotopies.
\item
$\eps^{(k,n)}_X$ restricts to a morphism
$e^{(k,n)}:Z^{(k)}\to Z^{(n)}$ for every $k\leq n$.
\item
$s.\trunc_ke^{(k,n)}:s.\trunc_kZ^{(n)}\to s.\trunc_kZ^{(k)}$
is an isomorphism, and $s.\trunc_kh^{(n)}=s.\trunc_kh^{(k)}$
for every $k\leq n$.
\item
The induced map
$Z^{(n)}[n]\to X[n]/\!\overset{p}{\sim}_{\partial\bDelta_n}$
is injective.
\end{enumerate}
\end{claim}
\begin{pfclaim} We construct inductively the sought
simplicial sets and morphisms, as follows. First,
set $Z^{(-1)}:=s.\emptyset$. Next, let $n\in\N$ be any
integer such that either $n=0$ or else
$Z^{(n-1)}\subset\Sk_{n-1}X$, $r^{(n-1)}$ and $h^{(n-1)}$
have already been defined.
Since $Z^{(n-1)}[n]$ consists entirely of degenerate
simplices (see remark \ref{rem_degenerate}(iv)), lemma
\ref{lem_degenerate}(i) implies that the natural projection
$Z^{(n-1)}[n]\to X[n]/\!\overset{p}{\sim}_{\partial\bDelta_n}$
is injective.

Now, for $n=0$, set $C_0:=X[0]/\!\overset{p}{\sim}_{\partial\bDelta_n}$,
and if $n>0$, let $C_n\subset X[n]/\!\overset{p}{\sim}_{\partial\bDelta_n}$
be the subset of all equivalence classes of elements
$x\in X[n]$ such that $\partial_ix\in Z^{(n-1)}[n-1]$ for
every $i=0,\dots,n$ (notice that $\partial_ix$ depends
only on the class of $x$ in $X[n]/\!\overset{p}{\sim}_{\partial\bDelta_n}$).
Pick a representative $x\in X[n]$ for every equivalence class
$\bar x\in C_n$ that does not lie in the image of $Z^{(n-1)}[n]$,
and let $R_n\subset X[n]$ be the subset of all such chosen
representatives; we denote by $T^{(n)}\subset s.\trunc_n X$
the $n$-truncated simplicial set such that
$$
T^{(n)}[n]=Z^{(n-1)}[n]\cup R_n
\qquad\text{and}\qquad
T^{(n)}[k]=Z^{(n-1)}[k]
\qquad
\text{for every $k<n$}.
$$
The inclusion map $T^{(n)}\to s.\trunc_n X$ induces a
morphism $\sk_nT^{(n)}\to\Sk_nX$, so we may set
$$
Z^{(n)}:=\Img\,(\sk_nT^{(n)}\to\Sk_nX)
$$
and we let $j^{(n)}:Z^{(n)}\to\Sk_nX$ be the inclusion map.
With this definition, we already see that conditions (d)
and (f) are fulfilled, as well the part of condition (e)
concerning $Z^{(n)}$. Next, we define the morphism
$r^{(n)}:\Sk_nX\to Z^{(n)}$. Namely, we set
$r^{(n)}[k]:=r^{(n-1)}[k]$ for every $k<n$, and the map
$r^{(n)}[n]$ is constructed as follows. If $n=0$, we let
$r^{(0)}[0]:X[0]\to Z^{(0)}[0]=R_0$ be the map that sends
any $x\in X[0]$ to the unique representative of the
equivalence class of $x$ in
$X[0]/\!\!\overset{p}{\sim}_{\partial\bDelta_n}$ that
lies in $R_0$. In case $n>0$, let $x\in\Sk_nX[n]=X[n]$
be any element; if $x\in Z^{(n)}[n]$, we let
$r^{(n)}[n](x):=x$, and if $x\notin Z^{(n)}[n]$ we let
$h_{\partial x}:\partial\bDelta_n\times\bDelta_1\to X$ be
the composition
$$
\partial\bDelta_n\times\bDelta_1
\xrightarrow{\ (\partial x)\times\one_{\bDelta_1}\ }
\Sk_{n-1}X\times\bDelta_1\xrightarrow{\ h^{(n-1)}\ }X.
$$
The pair $(h_{\partial x},x)$ amounts to a morphism
$\gamma_x$ fitting into a commutative diagram
$$
\xymatrix{
(\bDelta_n\times\bLambda^1_j)\cup
((\partial\bDelta_n)\times\bDelta_1) \ar[rrr]^-{\gamma_x} \ar[d]
& & & X \ar[d]^p \\
\bDelta_n\times\bDelta_1 \ar[r] &
\bDelta_n \ar[rr]^-{p\circ x} & & Y
}$$
whose horizontal (resp. vertical) unmarked arrow is the
natural projection (resp. the natural inclusion map).
Since $p$ is a fibration, proposition \ref{prop_saturate}
and lemma \ref{lem_satur-RLP}(i) imply that $\gamma_x$
extends to a morphism $h_x:\bDelta_n\times\bDelta_1\to X$
which, by inspection, is a $p$-fibrewise homotopy from $x$ to
$$
x':=h_x\circ(\one_{\bDelta_n}\times\eps_1):\bDelta_n\to X
$$
(notation of example \ref{ex_simplicial-sets}(iv)). Moreover,
$\partial_ix'\in Z^{(n-1)}[n-1]$ for every $i=0,\dots,n$, so
there exists a unique $x''\in T^{(n)}[n]$ such that
$x'\overset{p}{\sim}_{\partial\bDelta_n}x''$, and we set
$$
r^{(n)}[n](x):=x''.
$$
We have to check that the resulting system
$r^{(n)}[\bullet]:=(r^{(n)}[k]~|~k=0,\dots,n)$ is a morphism
$$
s.\trunc_nX\to T^{(n)}
$$
of $n$-truncated simplicial sets. If $n=0$, there is
nothing to show. If $n>0$, since the subsystem
$(r^{(n)}[k]~|~k=0,\dots,n-1)$ is already, by assumption,
a morphism $s.\trunc_{n-1}X\to s.\trunc_{n-1}Z^{(n-1)}$,
we have only to show that
$$
\partial(r^{(n)}[n](x))=r^{(n)}[n-1]\circ\partial x
\qquad
\text{for every $x\in X[n]$}.
$$
If $x\in T^{(n)}[n]$, the identity is obvious, since the
face and degeneracies of $T^{(n)}$ are the restrictions
of those of $X$. Hence, suppose $x\notin T^{(n)}[n]$;
then, by construction, $\partial(r^{(n)}[n](x))$ is
the morphism
$h_{\partial x}\circ(\one_{\partial\bDelta_n}\times\eps_1):
\partial\bDelta_n\to X$. The latter is the same as
$$
(\eta^{(n-1)}_X\circ j^{(n-1)}\circ r^{(n-1)})[n-1]\circ\partial x
$$
(because $h^{(n-1)}$ is a homotopy from $\eta_X^{(n-1)}$
to $\eta_X^{(n-1)}\circ j^{(n-1)}\circ r^{(n-1)}$) whence
the contention.

By adjunction, the morphism $r^{(n)}[\bullet]$ yields
the required morphism
$$
r^{(n)}:\Sk_nX\to Z^{(n)}.
$$
Also, again by adjunction, the identities of (a)
can be checked on the respective $n$-truncations,
where they are clear (details left to the reader).

It remains to exhibit the homotopy
$h^{(n)}:\Sk_nX\times\bDelta_1\to X$. However,
to any given $x\in\Sk_nX[n]\setminus Z^{(n)}[n]$, we
have already attached a $p$-fibrewise homotopy $h_x$
from $x$ to $x'$, and we have
$x'\overset{p}{\sim}_{\partial\bDelta_n}r^{(n)}[n](x)$; by lemma
\ref{lem_fibrewise-equiv}(ii), it follows that there
exists a $p$-fibrewise homotopy
$h'_x:\bDelta_n\times\bDelta_1\to X$ from $x$ to
$r^{(n)}[n](x)$ such that
\set\begin{equation}\label{eq_dontchg-boundary}
h'_x\circ(i_n\times\one_{\bDelta_1})=
h_x\circ(i_n\times\one_{\bDelta_1})=h_{\partial x}
\end{equation}
where $i_n:\partial\bDelta_n\to\bDelta_n$ is the inclusion
map. If $x\in Z^{(n)}[n]$, we have $x=r^{(n)}[n](x)$, so we
may take $h'_x$ to be a constant homotopy. There results a
well defined map
$$
\mu:X[n]\to s.\cHom(\bDelta_1,X)[n]
\qquad
x\mapsto h'_x.
$$
On the other hand, by remark \ref{rem_simpl-hom}(iv),
the homotopy $h^{(n-1)}$ corresponds to a morphism of
simplicial sets
$$
\nu^{(n-1)}:\Sk_{n-1}X\to s.\cHom(\bDelta_1,X).
$$
We claim that the datum of $\mu$ and the maps
$(\nu^{(n-1)}[k]~|~k=0,\dots,n-1)$ amounts to a
morphism of $n$-truncated simplicial sets
\set\begin{equation}\label{eq_truncamus}
s.\trunc_nX\to s.\trunc_n(s.\cHom(\bDelta_1,X)).
\end{equation}
The assertion comes down to checking the identity
\set\begin{equation}\label{eq_good-identity}
\partial_i\circ\mu(x)=\nu^{(n-1)}[n-1]\circ\partial_i(x)
\qquad
\text{for every $x\in X[n]$ and every $i=0,\dots,n$}
\end{equation}
and recall that the $i$-th face operator on
$s.\cHom(\bDelta_1,X)[n]$
$$
\partial_i:\Hom_{s.\Set}(\bDelta_n\times\bDelta_1,X)\to
\Hom_{s.\Set}(\bDelta_{n-1}\times\bDelta_1,X)
$$
is given by the rule :
$\phi\mapsto\phi\circ(\eps_i\times\one_{\bDelta_1})$
for every morphism $\phi:\bDelta_n\times\bDelta_1\to X$.

Hence, if $x\in Z^{(n)}[n]$, the morphism $\partial_i(h'_x)$
is the constant homotopy from $\partial_ix$ to itself; since
$h^{(n-1)}$ fulfills condition (c), this is the same as
$\nu^{(n-1)}[n-1](\partial_ix)$, so in this case
\eqref{eq_good-identity} holds. In case $x\notin Z^{(n)}[n]$,
identity \eqref{eq_dontchg-boundary} says that 
$\partial(h'_x)=h_{\partial x}$. But a simple inspection shows
that $h_{\partial x}=\nu^{(n-1)}[n-1]\circ\partial x$, so
\eqref{eq_good-identity} holds also in this case, and we
get the sought morphism \eqref{eq_truncamus}; by adjunction,
the latter corresponds to a morphism
$$
h^{(n)}:\Sk_nX\times\bDelta_1\to X
$$
as required, and clearly condition (e) for $h^{(n)}$ is
fulfilled by construction. Moreover, a simple inspection
shows that \eqref{eq_truncamus} fits into a commutative
diagram
$$
\xymatrix{
& & & s.\trunc_nX \\
\ar[rrru]^{\one_{s.\trunc_nX}\ \ \ \ }
s.\trunc_nX \ar[rrrd]_{s.\trunc_n(j^{(n)}\circ r^{(n)})\ \ \ \ }
\ar[rrr] & & & s.\trunc_n(s.\cHom(\bDelta_1,X))
\ar[u]_{s.\trunc_n(s.\cHom(\eps_0,X))}
\ar[d]^{s.\trunc_n(s.\cHom(\eps_1,X))} \\
& & & s.\trunc_nX
}$$
where we have used the natural identification
$X\isom s.\cHom(\bDelta_0,X)$ for the target of
the morphisms $s.\cHom(\eps_i,X)$ ($i=0,1$).
By adjunction, this diagram translates as our
condition (b) for $h^{(n)}$. To conclude the proof
of the claim, it remains to check condition (c).
However, in order to prove that $p\circ h^{(n)}$
is a constant homotopy it suffices -- by adjunction --
to show that the morphism
$\nu^{(n)}:\Sk_nX\to s.\cHom(\bDelta_1,X)$ corresponding
to $h^{(n)}$ fits into a commutative diagram
$$
\cD\quad : \quad
{\diagram \Sk_nX \ar[rrr]^-{\nu^{(n)}} \ar[d]_{p^{(n)}} & & &
s.\cHom(\bDelta_1,X) \ar[d]^{s.\cHom(\bDelta_1,p)} \\
Y=s.\cHom(\bDelta_0,Y) \ar[rrr]^-{s.\cHom(\pi,Y)}
& & & s.\cHom(\bDelta_1,Y)
\enddiagram}
$$
where $\pi$ is the unique morphism $\bDelta_1\to\bDelta_0$.
Again by adjunction, it then suffices to check that
$s.\trunc_n\cD$ commutes; the latter assertion follows
easily by inspecting the construction of the morphism
\eqref{eq_truncamus}. We argue similarly to show that
$h_Z^{(n)}$ is a constant homotopy : the assertion comes
down to checking the commutativity of the diagram
$$
\cD_Z\quad :\quad
{\diagram Z^{(n)} \ar[rrr]^{\nu^{(n)}\circ j^{(n)}} \ddouble
& & & s.\cHom(\bDelta_1,X) \\
Z^{(n)}=s.\cHom(\bDelta_0,Z^{(n)})
\ar[rrr]^-{s.\cHom(\pi,Z^{(n)})} & & &
s.\cHom(\bDelta_1,Z^{(n)})
\ar[u]_{s.\cHom(\bDelta_1,\eta^{(n)}_X\circ j^{(n)})}
\enddiagram}$$
and since the counit of adjunction
$\sk_nT^{(n)}\to Z^{(n)}$ is an epimorphism, we are
reduced to showing that $s.\trunc_n\cD_Z$ commutes,
which again follows by inspecting the construction
of \eqref{eq_truncamus} : details left to the reader.
\end{pfclaim}

Now, if the system $(Z^{(n)},j^{(n)},r^{(n)},h^{(n)}~|~n\in\N)$
fulfills conditions (a)--(f) of claim \ref{cl_approx-def-retract},
we let
$$
Z:=\colim_{n\in\N}Z^{(n)}
$$
where the transition maps $Z^{(k)}\to Z^{(n)}$ are
the morphisms $e^{(k,n)}$, for every $k,n\in\N$ with
$k\leq n$. In view of \eqref{eq_isom-trunc}, the colimit
of the system of monomorphisms $(j^{(n)}~|~n\in\N)$ (resp.
epimorphisms $(r^{(n)}~|~n\in\N)$) is therefore a monomorphism
$$
j:Z\to X
\qquad
\text{(resp. an epimorphism $r:X\to Z$)}
$$
such that $r\circ j=\one_Z$ and $p\circ j\circ r=p$.
It is also clear that there exists a unique homotopy
$h:X\times\bDelta_1\to X$ from $\one_X$ to $j\circ r$,
such that
$$
s.\trunc_n(h)=s.\trunc_n(h^{(n)})
\qquad
\text{for every $n\in\N$}
$$
and then the resulting datum $(j,r,h)$ is the sought
strong deformation retract. It remains to check that
$p\circ j$ is a minimal fibration. First, it follows
from lemma \ref{lem_satur-RLP}(iv) that $p\circ j$ is
a fibration. Lastly, fix $n\in\N$, and suppose that
$z,z'\in Z[n]$ are any two elements such that
$z\overset{p\circ j}{\sim}_{\partial\bDelta_n}z'$; then
$j(z)\overset{p}{\sim}_{\partial\bDelta_n}j(z')$, so $z=z'$,
by condition (f) of claim \ref{cl_approx-def-retract}.
\end{proof}

\begin{proposition} 
Consider a commutative diagram of simplicial sets
$$
\xymatrix{ X' \ar[rr]^-f \ar[rd]_{p'} & & X \ar[ld]^p \\
& Y
}$$
such that $p$ and $p'$ are minimal fibrations, and $f$
is a fibrewise homotopy equivalence. Then $f$ is an
isomorphism of simplicial sets.
\end{proposition}
\begin{proof} The proposition is an immediate consequence
of the following more precise :

\begin{claim} In the situation of the proposition,
suppose that :
\begin{enumerate}
\alphaenu
\item
$f':X'\isom X$ is an isomorphism such that $p\circ f'=p'$
and $f\overset{p}\sim f'$.
\item
$p$ is a minimal fibration.
\end{enumerate}
Then $f$ is an isomorphism as well.
\end{claim}
\begin{pfclaim}[] We show first that $f[n]$ is injective
for every $n\in\N$. We argue by induction on $n\in\N$.
For $n=0$, notice that $\partial\bDelta_0=s.\emptyset$;
hence, if $a,b\in X'[0]$ and $f[0](a)=f[0](b)$, assumption
(a) yields :
$$
f'[0](a)\overset{p}{\sim}_{\partial\bDelta_0}f[0](a)
\qquad
f'[0](b)\overset{p}{\sim}_{\partial\bDelta_0}f[0](b)
$$
whence $f'[0](a)\overset{p}{\sim}_{\partial\bDelta_0}f'[0](b)$
by virtue of lemma \ref{lem_fibrewise-equiv}(i); then
assumption (b) implies that $f'[0](a)=f'[0](b)$, and
finally $a=b$, as stated. Next, suppose that $n>0$
and $f[n-1]$ is already known to be injective, and let
$a,b\in X'[n]$ be any two elements with $f[n](a)=f[n](b)$;
condition (a) says that there exists a $p$-fibrewise
homotopy $h_a$ (resp. $h_b$) : $\bDelta_n\times\bDelta_1\to X$
from $f[n](a)$ to $f'[n](a)$ (resp. from $f[n](b)$ to
$f'[n](b)$), and on the other hand, the inductive assumption
implies that $\partial_ia=\partial_ib$, for every $i=0,\dots,n$,
so the restrictions of $h_a$ and $h_b$ are both equal
to the same morphism $w:\partial\bDelta_n\times\bDelta_1\to X$.
There follows a commutative diagram
$$
\xymatrix{ (\bDelta_n\times\bLambda_2^2)\cup
(\partial\bDelta_n\times\bDelta_2) \ar[r]^-u \ar[d]
& X \ar[d]^p \\
\bDelta_n\times\bDelta_2 \ar[r]^-v & Y
}$$
such that
$$
u\circ(\one_{\bDelta_n}\times\eps_1)=h_a
\qquad
u\circ(\one_{\bDelta_n}\times\eps_0)=h_b
$$
and the restriction of $u$ to $\partial\bDelta_n\times\bDelta_2$
equals $w\circ(\one_{\partial\bDelta_n}\times\eta_0)$. Also,
$v$ is the composition of the projection
$\bDelta_n\times\bDelta_2\to\bDelta_n$ and the morphism
$(p\circ f)[n](a):\bDelta_n\to Y$. Then there exists
a morphism $t:\bDelta_n\times\bDelta_2\to X$ that extends
$u$ and lifts $v$; it is easily seen that
$t\circ(\bDelta_n\times\eps_2):\bDelta_n\times\bDelta_1\to X$
is a $p$-fibrewise homotopy from $f'[n](a)$ to $f'[n](b)$
relative to $\partial\bDelta_n$, and since $p$ is minimal,
we conclude that $f'[n](a)=f'[n](b)$, so finally $a=b$,
as required.

Next, we show by induction on $n\in\N$ that $f[k]$ is
surjective, for every $k<n$. If $n=0$, there is nothing
to show, hence suppose that $n>0$, and let $x\in X[n]$
be any element; by inductive assumption, for every
$i=0,\dots,n$ there exists $a_i\in X'[n-1]$ such that
$f[n-1](a_i)=\partial_ix$. Since the injectivity of
$f$ has already been established, we see easily that
$\partial_ia_j=\partial_{j-1}a_i$ whenever $0\leq i<j\leq n$,
so the system $(a_0,\dots,a_n)$ determines a unique
morphism
$$
w':\partial\bDelta_n\to X'
$$
such that $f\circ w'$ is the restriction of $x:\bDelta_n\to X$
to $\partial\bDelta_n$ (see example \ref{ex_simplicial-sets}(iv)).
Now, let $h:X'\times\bDelta_1\to X$ be a $p$-fibrewise homotopy
from $f$ to $f'$; we deduce a commutative diagram
$$
\xymatrix{
(\partial\bDelta_n\times\bDelta_1)\cup(\bDelta_n\times\bLambda_1^0)
\ar[r]^-{u'} \ar[d] & X \ar[d]^p \\
\bDelta_n\times\bDelta_1 \ar[r]^-{v'} & Y
}$$
such that the restriction of $u'$ to $\partial\bDelta_n\times\bDelta_1$
(resp. to $\bDelta_n\times\bLambda_1^0$) agrees with
$h\circ(w'\times\one_{\bDelta_1})$ (resp. with the composition
of the projection $\bDelta_n\times\bLambda_1^0\isom\bDelta_n$
and the morphism $x$). Also, $v'$ is a constant homotopy.
There follows a morphism $t':\bDelta_n\times\bDelta_1\to X$
that extends $u'$ and lifts $v'$. Let
$$
x':\bDelta_n\to X'
$$
be the unique morphism such that
$f'\circ x'=t'\circ(\one_{\bDelta_n}\times\eps_1)$.
By construction, we have
$$
f'\circ x'\circ i_n=f'\circ w'
$$
(where $i_n:\partial\bDelta_n\to\bDelta_n$ is the
natural inclusion); whence $x'\circ i_n=w'$, and we
obtain yet another commutative diagram
$$
\xymatrix{
(\bDelta_n\times\bLambda_2^0)\cup(\partial\bDelta_n\times\bDelta_2)
\ar[r]^-{u''} \ar[d] & X \ar[d]^p \\
\bDelta_n\times\bDelta_2 \ar[r]^-{v''} & Y
}$$
such that
$$
u''\circ(\one_{\bDelta_n}\times\eps_1)=t'
\qquad
u''\circ(\one_{\bDelta_n}\times\eps_2)=
h\circ(x'\times\one_{\bDelta_1})
$$
and where $v''$ factors through the projection
$\bDelta_n\times\bDelta_2\to\bDelta_n$. As usual,
we find a morphism $t'':\bDelta_n\times\bDelta_2\to X$
that lifts $v''$ and extends $u''$, and it is easily
seen that $t''\circ(\bDelta_n\times\eps_0)$ is a
$p$-fibrewise homotopy from $x$ to $f[n](x')$ relative
to $\partial\bDelta_n$. By the minimality of $p$,
it follows that $x=f[n](x')$, whence the sought
surjectivity of $f[n]$.
\end{pfclaim}
\end{proof}

\begin{definition}
Let $f:X\to Y$ be any morphism of fibrant simplicial
sets. We say that $f$ is a {\em weak equivalence}, if
the induced map $\pi_n(f,\xi)$ is an isomorphism,
for every $\xi\in X[0]$ and every $n\in\N$.
\end{definition}

\begin{proposition}\label{prop_weak-equiv}
Let $f:X\to Y$ and $g:Y\to Z$ be any two morphisms
of fibrant simplicial sets. The following holds :
\begin{enumerate}
\item
If $f$ has the right lifting property with respect to the
natural monomorphism $i_n:\partial\bDelta_n\to\bDelta_n$
for every $n\in\N$, then $f$ is both a fibration and a weak
equivalence.
\item
If $f$ is a homotopy equivalence, then it is a weak
equivalence.
\item
If any two of the three morphisms $f,g,g\circ f$ is a
weak equivalence, then so is the third.
\end{enumerate}
\end{proposition}
\begin{proof}(i): It follows easily from lemma
\ref{lem_satur-RLP}(i) and remark \ref{rem_skeleta}
that $f$ has the right lifting property with respect
to every monomorphism of simplicial sets (details
left to the reader). Especially, $f$ is a fibration.
It also follows easily that $f[0]$ is surjective.
Moreover, suppose that $a,b\in X[0]$ and $[f(a)]=[f(b)]$
in $\pi_0(Y)$, and pick a homotopy $h:\bDelta_1\to Y$
from $f(a)$ to $f(b)$; there follows a commutative
diagram
$$
\xymatrix{
\partial\bDelta_1 \ar[d]_{i_1} \ar[r]^-c & X \ar[d]^f \\
\bDelta_1 \ar[r]^-h & Y
}$$
such that the restriction of $c$ to $\bLambda^0_1$ (resp.
to $\bLambda^1_1$) agrees with $a$ (resp. with $b$).
By assumption, $c$ extends to a morphism $h':\bDelta_1\to X$
such that $f\circ h'=h$. We conclude that $\pi_0(f)$
is bijective. Next, let $\xi\in X[0]$ be any vertex,
and set $F:=f^{-1}(f(\xi))$. By lemma \ref{lem_satur-RLP}(ii),
the unique morphism $F\to\bDelta_0$ has the right lifting
properties with respect to $i_n$, for every $n\in\N$.
By the foregoing, we know already that $\pi_0(F)$ is
a set of cardinality one; moreover, we have :

\begin{claim}\label{cl_F-is-trivial}
$\pi_n(F,\xi)=0$ for every $n\geq 1$.
\end{claim}
\begin{pfclaim} Let $a\in F(\xi,n)$ be any element; then
there exists a unique morphism $b:\partial\bDelta_{n+1}\to F$
whose restriction to $\bLambda^0_{n+1}$ factors through
$\xi$, and whose composition with
$\eps_0:\bDelta_n\to\partial\bDelta_{n+1}$ agrees with $a$.
By assumption, $b$ extends to a morphism $c:\bDelta_{n+1}\to X$,
and then the claim follows from lemma \ref{lem_hom-groups}.
\end{pfclaim}

The proposition now follows easily from theorem
\ref{th_homotopy-seq-fib} and claim \ref{cl_F-is-trivial}.

(ii): Let $f':Y\to X$ be a morphism and $h:X\times\bDelta_1\to X$
a homotopy from $\gamma_0:=\one_X$ to $\gamma_1:=f'\circ f$.
By adjunction, $h$ corresponds to a morphism
$\beta:X\to s.\cHom(\bDelta_1,X)$ fitting into the commutative
diagrams
$$
{\diagram
X \ar[r]^-\beta \ar[rd]_{\gamma_i} &
s.\cHom(\bDelta_1,X) \ar[d]^{s.\cHom(\eps_i,X)} \\
& s.\cHom(\bLambda^i_1,X)=X
\enddiagram}
\qquad
i=0,1.
$$
Recall that $s.\cHom(\bDelta_1,X)$ is fibrant (corollary
\ref{cor_fibrant-hom}(i)).

\begin{claim}\label{cl_these-are-weak}
$\beta$ and $\gamma_1$ are weak equivalences.
\end{claim}
\begin{pfclaim} As for $\beta$, since $\gamma_0$ is an
isomorphism, it suffices to show that $p_0:=s.\cHom(\eps_0,X)$
is a weak equivalence. To this aim, by (i), we are reduced
to checking that $p_0$ has the right lifting property
with respect to the monomorphisms $i_n$, for every $n\in\N$.
Thus, let $a:\partial\bDelta_n\to s.\cHom(\bDelta_1,X)$
be a morphism such that $p_0\circ a$ extends to a
morphism $b:\bDelta_n\to X$; by adjunction, $a$ and
$b$ determine a unique morphism
$(\bDelta_n\times\bLambda^0_1)\cup
(\partial\bDelta_n\times\bDelta_1)\to X$ which, since
$X$ is fibrant, extends to a morphism
$\bDelta_n\times\bDelta_1\to X$. The latter corresponds,
by adjunction, to a unique morphism
$\bDelta_n\to s.\cHom(\bDelta_1,X)$ which is the sought
extension of $a$.

The same argument shows as well that $s.\cHom(\eps_1,X)$
is a weak equivalence, and so the same follows for
$\gamma_1$, as stated.
\end{pfclaim}

Now, if $f$ is a homotopy equivalence, we have also a
homotopy $h':Y\times\bDelta_1\to Y$ from $\one_Y$ to
$f\circ f'$. Taking into account claim \ref{cl_these-are-weak},
we deduce that both
$$
\pi_n(g,f(\xi))\circ\pi_n(f,\xi)
\quad\text{and}\quad
\pi_n(f,f'(\xi'))\circ\pi_n(f',\xi')
$$
are bijections, for every $n\in\N$, every $\xi\in X[0]$
and every $\xi'\in Y[0]$. The assertion follows.

(iii): The only non-trivial assertion is that $g$
is a weak equivalence if both $f$ and $g\circ f$ are.
However, in this case it is clear that $\pi_0(g)$
is an isomorphism, so it remains only to check that
$\pi_n(g,y)$ is an isomorphism for every $y\in Y[0]$.
Thus, fix such vertex $y$; since $\pi_0(f)$ is a
bijection, there exist $x\in X[0]$ and a morphism
$[\alpha]:y\to y':=f[0](x)$ in $\underline\pi(Y)$
(notation of \eqref{subsec_device-groupoid}).
Set $z:=g[0](y)$, $z':=g[0](z')$ and $\beta:=g[1](\alpha)$;
there follows a commutative diagram (see remark
\ref{rem_sothat})
$$
\xymatrix{
\pi_n(Y,y) \ar[rr]^-{\pi_n(Y,[\alpha])} \ar[d]_{\pi_n(g,y)} & &
\pi_n(Y,y') \ar[d]^{\pi_n(g,y')} & & &
\pi_n(X,x) \ar[lll]_-{\pi_n(f,x)} \ar[llld]^-{\pi_n(g\circ f,x)} \\
\pi_n(Z,z) \ar[rr]^-{\pi_n(Z,[\beta])} & & \pi_n(Z,z')
}$$
from which we see first that $\pi_n(g,y')$ is an isomorphism,
and then the assertion follows.
\end{proof}

The following is the main result of this section :

\begin{theorem}[Whitehead] Every weak equivalence of
fibrant simplicial sets is a homotopy equivalence.
\end{theorem}
\begin{proof} We begin with the following special case :

\begin{claim}\label{cl_minimal-case}
Let $f:X\to Y$ be any minimal fibration of fibrant
simplicial sets, and suppose that $f$ is a weak
equivalence. Then $f$ is an isomorphism.
\end{claim}
\begin{pfclaim} We consider first the case where $Y=\bDelta_0$
(notice that $\bDelta_0$ is trivially a fibrant simplicial set),
so that $f$ is the unique morphism $X\to\bDelta_0$. We show by
induction that $X[n]$ contains exactly one element, for every
$n\in\N$. Say first that $n=0$; since $\pi_0(\bDelta_0)$ is
non-empty, we may find a vertex $\xi\in X[0]$, and we set
$\xi':=f[0](\xi)$. Let $a\in X[0]$ be any other vertex; since
$f[0](a)=\xi'$, by assumption we may find a homotopy from
$a$ to $\xi$, and therefore $a=\xi$, since $X$ is minimal. Next,
suppose that $n>0$, and the assertion is already known for every
integer $<n$; let $b\in X[n]$ be any simplex, and notice
that $\partial b$ must factor through $\xi$, by inductive
assumption, so $b\in X(\xi,n)$. However, obviously
$\pi_n(\bDelta_0,\xi')=0$, hence $b\sim_{\partial\bDelta_n}\xi$,
and finally $b=\xi$, again by the minimality of $X$.

Next, consider an arbitrary minimal fibration $f$, fix
$\xi'\in Y[0]$ and pick any $\xi\in f[0]^{-1}(\xi')$ (notice
that $Y[0]$ and $f[0]^{-1}(\xi')$ are non-empty, since
$Y$ is fibrant and $f$ is a fibration). Set $F:=f^{-1}(\xi')$;
theorem \ref{th_homotopy-seq-fib}(i) implies that
$\pi_n(F,\xi)$ has cardinality equal to one, for every
$n\in\N$. On the other hand, since $f$ is minimal, it is
easily seen that the same holds for $F$ (details left to
the reader) therefore $F\simeq\bDelta_0$, by the foregoing
case. Since $\xi'$ is arbitrary, the claim follows.
\end{pfclaim}

Now, let $f:X\to Y$ be any weak equivalence, with $X$ and
$Y$ fibrant simplicial sets. By theorem \ref{th_E-infinity},
there exist an anodyne extension $i:X\to E$ and a
fibration $p:E\to Y$ such that $p\circ i=f$. Since
$Y$ is fibrant, the same holds for $E$.

\begin{claim}\label{cl_i-is-retract}
The morphism $i$ is a strong deformation retract of $E$.
\end{claim}
\begin{pfclaim} Since $i$ is anodyne and $X$ is fibrant,
the identity $\one_X$ extends to a morphism $r:E\to X$.
Next, consider the unique morphism
$$
t:(X\times\bDelta_1)\cup(E\times\partial\bDelta_1)\to E
$$
whose restriction to $X\times\bDelta_1$ is the composition
of $i$ with the projection onto $X$, and whose restriction
with $E\times\bLambda^0_1$ (resp. with $E\times\bLambda^1_1$)
agrees with $\one_E$ (resp. with $i\circ r$). Since $E$
is fibrant, $t$ extends to a morphism $h:E\times\bDelta_1\to E$
(corollary \ref{cor_saturate}), and it is easily seen that
$(i,r,h)$ is a strong deformation retract for the unique
morphism $E\to\bDelta_0$.
\end{pfclaim}

By claim \ref{cl_i-is-retract} and remark
\ref{rem_equivalences}(ii) we see that $i$ is a homotopy
equivalence, hence it suffices to check that the same
holds for $p$. Moreover, $i$ is a weak equivalence
(proposition \ref{prop_weak-equiv}(ii)), hence the same
holds for $p$ (proposition \ref{prop_weak-equiv}(iii)).
Hence, we may replace $X$ by $E$ and $f$ by $p$, and
assume from start that $f$ is also a fibration. We may
then find a strong deformation retract $(j:F\to X, r,h)$
of $f$ such that $q:=f\circ j$ is a minimal fibration.
Arguing as in the foregoing, we see that $j$ is a weak
equivalence, and therefore the same holds for $q$
(again, by proposition \ref{prop_weak-equiv}(iii));
then claim \ref{cl_minimal-case} says that $q$ is an
isomorphism. Lastly, $r$ is a homotopy equivalence
(remark \ref{rem_equivalences}(ii)), so the same holds
for $f=q\circ r$.
\end{proof}

\subsection{Graded rings}
To motivate the results of this section, let us review briefly
the well known correspondance between $\Gamma$-gradings on a
module $M$ (for a given commutative group $\Gamma$), and
actions of the diagonalizable group $D(\Gamma)$ on $M$.

\sset\subsubsection{Graded rings}\label{subsec_monoid-graded}
Let $S$ be a scheme; on the category $\Sch/S$ of $S$-schemes
we have the presheaf of rings :
$$
\cO_{\!\Sch/S}:\Sch/S\to\Z\Alg
\qquad
(X\to S)\mapsto\Gamma(X,\cO_{\!X}).
$$
Also, let $M$ be an $\cO_{\!S}$-module; following
\cite[Exp.I, D\'ef.4.6.1]{SGA3}, we attach to $M$ the
$\cO_{\!\Sch/S}$-module $\cW_M$ given by the rule :
$(f:X\to S)\mapsto\Gamma(X,f^*M)$. If $M$ and $N$ are two
quasi-coherent $\cO_{\!S}$-modules, and $f:S'\to S$ an affine
morphism, it is easily seen that
\set\begin{equation}\label{eq_many-ws}
\Gamma(S',\cHom_{\cO_{\Sch/S}}(\cW_M,\cW_N))=
\Hom_{\cO_{\!S}}(M,f_*\cO_{\!S'}\otimes_{\cO_{\!S}}N)
\end{equation}
(see \cite[Exp.I, Prop.4.6.4]{SGA3} for the details).

Let $f:G\to S$ be a group $S$-scheme, {\em i.e.} a group object
in the category $\Sch/S$. If $f$ is affine, we say that $G$
is an {\em affine group $S$-scheme}; in that case, the
mutiplication law $G\times_SG\to G$ and the unit section
$S\to G$ correspond respectively to morphisms of $\cO_{\!S}$-algebras
$$
\Delta_G:f_*\cO_{\!G}\to f_*\cO_{\!G}\otimes_{\cO_{\!S}}f_*\cO_{\!G}
\qquad
\eps_G:\cO_{\!S}\to f_*\cO_{\!G}
$$
which make commute the diagram :
$$
\xymatrix{ f_*\cO_{\!G} \ar[rrr]^-{\Delta_G} \ar[d]_{\Delta_G} & & &
f_*\cO_{\!G}\otimes_{\cO_{\!S}}f_*\cO_{\!G}
\ar[d]^{\one_{f_*\cO_{\!S}}\otimes\Delta_G} \\
f_*\cO_{\!G}\otimes_{\cO_{\!S}}f_*\cO_{\!G}
\ar[rrr]^-{\Delta_G\otimes\one_{f_*\cO_{\!S}}} & & &
f_*\cO_{\!G}\otimes_{\cO_{\!S}}f_*\cO_{\!G}\otimes_{\cO_{\!S}}f_*\cO_{\!G}
}$$
as well as a similar diagram, which expresses the unit property
of $\eps_G$ : see \cite[Exp.I, \S4.2]{SGA3}.

\begin{example}
Let $G$ be any commutative group. The presheaf of groups
$$
D_S(G):\Sch/S\to\Z\Mod
\qquad
(X\to S)\mapsto\Hom_{\Z\Mod}(G,\cO^\times_{\!X}(X))
$$
is representable by an affine group $S$-scheme $D_S(G)$,
called the {\em diagonalizable group scheme} attached to $G$.
Explicitly, if $S=\Spec\,R$ is an affine scheme, the underlying
$S$-scheme of $D_S(G)$ is $\Spec\,R[G]$, and the group law is given
by the map of $R$-algebras
$$
\Delta_G:R[G]\to R[G]\otimes_RR[G]\isom R[G\times G]
\qquad
g\mapsto(g,g)
\quad
\text{for every $g\in G$}
$$
with unit $\eps_G:R[G]\to R$ given by the standard augmentation
(see \cite[Exp.I, \S4.4]{SGA3}). For a general scheme $S$,
we have $D_S(G)=D_{\Spec\,\Z}(G)\times_{\Spec\,\Z}S$ (with the
induced group law and unit section).
\end{example}

\begin{definition} Let $M$ be an $\cO_{\!S}$-module, $G$ a
group $S$-scheme. A {\em $G$-module structure} on
$M$ is the datum of a morphism of presheaves of groups on
$\Sch/S$ :
$$
h_G\to{\cA\!ut}_{\cO_{\!\Sch/S}}(\cW_M)
$$
(where $h_G$ denotes the Yoneda embedding : see
\eqref{subsec_yoneda}).
\end{definition}

\sset\subsubsection{}\label{subsec_G-module-affine}
Suppose now that $f:G\to S$ is an affine group $S$-scheme,
and $M$ a quasi-coherent $\cO_{\!S}$-module; in view of
\eqref{eq_many-ws}, a $G$-module structure on $M$ is then
the same as a map of $\cO_{\!S}$-modules
$$
\mu_M:M\to f_*\cO_{\!G}\otimes_{\cO_{\!S}}M
$$
which makes commute the diagrams :
$$
\xymatrix{ M \ar[rrr]^-{\mu_M} \ar[d]_{\mu_M} & & &
f_*\cO_{\!G}\otimes_{\cO_{\!S}}M \ar[d]^{\Delta_G\otimes\one_M} &
M \ar[rr]^-{\mu_M} \ar[rrd]_{\one_M} & &
f_*\cO_{\!G}\otimes_{\cO_{\!S}}M \ar[d]^{\eps_G\otimes\one_M} \\
f_*\cO_{\!G}\otimes_{\cO_{\!S}}M
\ar[rrr]^-{\one_{f_*\cO_{\!G}}\otimes\mu_M} & & &
f_*\cO_{\!G}\otimes_{\cO_{\!S}}f_*\cO_{\!G}\otimes_{\cO_{\!S}}M &
& & M.
}$$

\begin{example}\label{ex_diagonal-module}
Let $\Gamma$ be a commutative group; a $D_S(\Gamma)$-module
structure on a quasi-coherent $\cO_{\!S}$-module $M$ is the datum
of a morphism of $\cO_{\!S}$-modules
$$
\mu_M:M\to\cO_{\!S}[\Gamma]\otimes_{\cO_{\!S}}M=\Z[\Gamma]\otimes_\Z M
$$
which makes commute the diagrams of \eqref{subsec_G-module-affine}.
If $S=\Spec\,R$ is affine, $M$ is associated with an $R$-module
which we denote also by $M$; in this case, $\mu_M$ is the same
as a system $(\mu_M^{(\gamma)}~|~\gamma\in\Gamma)$ of $\cO_{\!S}$-linear
endomorphisms of $M$, such that :

\begin{itemize}
\item
for every $x\in M$, the subset
$\{\gamma\in\Gamma~|~\mu_M^{(\gamma)}(x)\neq 0\}$ is finite.
\item
$\mu_M^{(\gamma)}\circ\mu_M^{(\tau)}=\delta_{\gamma,\tau}\cdot\one_M$
\ \ and\ \ $\sum_{\gamma\in\Gamma}\mu_M=\one_M$.
\end{itemize}
In other words, the $\mu_M^{(\gamma)}$ form an orthogonal system
of projectors of $M$, summing up to the identity $\one_M$.
This is the same as the datum of a $\Gamma$-grading on $M$ : namely,
for a given $D_S(\Gamma)$-module structure $\mu_M$, one lets
$$
\gr_\gamma M:=\mu_M^{(\gamma)}(M)
\qquad
\text{for every $\gamma\in\Gamma$}
$$
and conversely, given a $\Gamma$-grading $\gr_\bullet M$ on $M$,
one defines $\mu_M$ as the $R$-linear map given by the rule :
$x\mapsto\gamma\otimes x$ for every $\gamma\in\Gamma$ and every
$x\in\gr_\gamma M$.
\end{example}

\sset\subsubsection{}\label{subsec_corres-diagonal}
Suppose now that $g:X\to S$ is an affine $S$-scheme, and $f:G\to S$
an affine group $S$-scheme. A {\em $G$-action\/} on $X$ is
a morphism of presheaves of groups :
$$
h_G\to\cA\!ut_{\Sch/S^\wedge}(h_X).
$$
(notation of \eqref{subsec_yoneda}); the latter is the
same as a morphism of $S$-schemes
\set\begin{equation}\label{eq_no-nam-action}
G\times_SX\to X
\end{equation}
inducing a $G$-module structure on $g_*\cO_{\!X}$ :
$$
g_*\cO_{\!X}\to f_*\cO_{\!G}\otimes_{\cO_{\!S}}g_*\cO_{\!X}
$$
which is also a morphism of $\cO_{\!S}$-algebras. For instance,
if $S=\Spec\,R$ is affine, and $G=D_S(\Gamma)$ for an abelian
group $\Gamma$, we may write $X=\Spec\,A$ for some $R$-algebra
$B$, and in view of example \ref{ex_diagonal-module}, the $G$-action
on $X$ is the same as the datum of a $\Gamma$-graded $R$-algebra
structure on $B$, in the sense of the following :

\begin{definition}\label{def_Gamma-graded-algs}
Let $(\Gamma,+)$ be a commutative monoid, $R$ a ring.

(i)\ \
A {\em $\Gamma$-graded (associative, unital) $R$-algebra\/}
is a pair $\underline B:=(B,\gr_\bullet B)$ consisting of an
associative unital $R$-algebra $B$ and a $\Gamma$-grading
$B=\bigoplus_{\gamma\in\Gamma}\gr_\gamma B$ of the $R$-module $B$,
such that
$$
\gr_\gamma B\cdot\gr_{\gamma'}B\subset\gr_{\gamma+\gamma'}B
\qquad
\text{for every $\gamma,\gamma'\in\Gamma$}.
$$
A morphism of $\Gamma$-graded $R$-algebras is a map of $R$-algebras
which is compatible with the gradings, in the obvious way.

(ii)\ \
Let $\underline B:=(B,\gr_\bullet B)$ be a $\Gamma$-graded
$R$-algebra. A {\em $\Gamma$-graded (left) $\underline B$-module\/}
is a datum $\underline M:=(M,\gr_\bullet M)$ consisting of
a $B$-module $M$ and a $\Gamma$-grading
$M=\bigoplus_{\gamma\in\Gamma}\gr_\gamma M$ of the $R$-module
underlying $M$, such that
$$
\gr_\gamma B\cdot\gr_{\gamma'}M\subset\gr_{\gamma+\gamma'}M
\qquad
\text{for every $\gamma,\gamma'\in\Gamma$}.
$$
A morphism of $\Gamma$-graded $\underline B$-module is a
map $f:M\to N$ of $B$-modules, such that
$f(\gr_\gamma M)\subset\gr_\gamma N$, for every $\gamma\in\Gamma$.
Likewise we define
{\em $\Gamma$-graded right $\underline B$-modules\/} and
{\em $\Gamma$-graded $\underline B$-bimodules}, and their
respective morphisms.

(iii)\ \
A {\em graded ideal} of $\underline B$ is a $\Gamma$-graded
sub-$\underline B$-bimodule $(I,\gr_\bullet I)$ of $\underline B$.
Obviously, in this case we must have $\gr_\gamma I=I\cap\gr_\gamma B$
for every $\gamma\in\Gamma$.

(iv)\ \
If $f:\Gamma'\to\Gamma$ is any morphism of commutative monoids,
and $\underline M$ is any $\Gamma$-graded $R$-module, we define
the $\Gamma'$-graded $R$-module $\Gamma'\times_\Gamma\underline M$
by setting
$$
\gr_\gamma(\Gamma'\times_\Gamma M):=\gr_{f(\gamma)}M
\qquad
\text{for every $\gamma\in\Gamma'$}.
$$
Notice that if $\underline B:=(B,\gr_\bullet B)$ is a
$\Gamma$-graded $R$-algebra, then $\Gamma'\times_\Gamma\underline B$
is a $\Gamma'$-graded $R$-algebra, with multiplication
and addition laws induced by those of $B$, in the obvious
way. Likewise, if $(M,\gr_\bullet M)$ is a $\Gamma$-graded
$\underline B$-module, then $\Gamma'\times_\Gamma M$ is naturally
a $\Gamma'$-graded $\Gamma'\times_\Gamma\underline B$-module.

(v)\ \
Furthermore, if $N$ is any $\Gamma'$-graded $R$-module, we define
the $\Gamma$-graded $R$-module $N_{/\Gamma}$ whose underlying
$R$-module is the same as $N$, and whose grading is given by the
rule
$$
\gr_\gamma(N_{/\Gamma}):=
\bigoplus_{\gamma'\in f^{-1}(\gamma)}\gr_\gamma N.
$$
Just as in (iii), if $\underline C:=(C,\gr_\bullet C)$ is a
$\Gamma'$-graded $R$-algebra, then we get a $\Gamma$-graded
$R$-algebra $\underline C_{/\Gamma}$, whose underlying $R$-algebra
is the same as $C$. Lastly, if $(N,\gr_\bullet N)$ is a
$\Gamma'$-graded $\underline C$-module, then $N_{/\Gamma}$
is a $\Gamma$-graded $\underline C_{/\Gamma}$-module.
\end{definition}

In this section we will consider only {\em commutative}
graded algebras, but in the next section
\eqref{subsec_diff-grad-alg} we will encounter certain
special classes of non-commutative graded algebras as well.

\begin{example}\label{ex_actions-are-thetas}
(i)\ \
For instance, the $R$-algebra $R[\Gamma]$ is naturally a
$\Gamma$-graded $R$-algebra, when endowed with the $\Gamma$-grading
such that $\gr_\gamma R[\Gamma]:=\gamma R$ for every $\gamma\in\Gamma$.

(ii)\ \
Suppose that $\Gamma$ is an integral monoid. Then, to a
$\Gamma$-graded $R$-algebra $B$, the correspondance described in
\eqref{subsec_corres-diagonal} attaches a $D_S(\Gamma^\gp)$-action
on $\Spec\,B$ (where $S:=\Spec\,R$), given by the map of $R$-algebras
$$
\theta_B:B\to B[\Gamma]\subset B[\Gamma^\gp]
\quad :\quad
b\mapsto b\cdot\gamma
\qquad
\text{for every $\gamma\in\Gamma$, and every $b\in\gr_\gamma B$.}
$$
\end{example}

\begin{remark}\label{rem_etale-and-smooth}
(i)\ \
Let $R$ be a ring; consider a cartesian diagram of monoids
$$
\xymatrix{ \Gamma_3 \ar[r] \ar[d] & \Gamma_1 \ar[d] \\
           \Gamma_2 \ar[r] & \Gamma_0
}$$
and let $B$ be any $\Gamma_1$-graded $R$-algebra. A simple
inspection of the definitions yields an identity of
$\Gamma_2$-graded $R$-algebras :
$$
\Gamma_2\times_{\Gamma_0}B_{/\Gamma_0}=
(\Gamma_3\times_{\Gamma_1}B)_{/\Gamma_2}.
$$

(ii)\ \
Suppose that $\Gamma$ is a finite abelian group,
whose order is invertible in $\cO_{\!S}$. Then $D_S(\Gamma)$
is an \'etale $S$-scheme. Indeed, in light of
\eqref{eq_push-out-tensor}, the assertion is reduced to the
case where $\Gamma=\Z/n\Z$ for some integer $n>0$ which is
invertible in $\cO_{\!S}$. However, $R[\Z/n\Z]\simeq R[T]/(T^n-1)$,
which is an \'etale $R$-algebra, if $n\in R^\times$.

(iii)\ \
More generally, suppose that $\Gamma$ is a finitely generated
abelian group, such that the order of its torsion subgroup
is invertible in $\cO_{\!S}$. Then we may write
$\Gamma=L\oplus\Gamma_{\!\mathrm{tor}}$, where $L$ is a free
abelian group of finite rank, and $\Gamma_{\!\mathrm{tor}}$
is a finite abelian group as in (ii). In view of
\eqref{eq_push-out-tensor} and (ii), we conclude that
$D_S(\Gamma)$ is a smooth $S$-scheme in this case.

(iv)\ \
Let $\Gamma$ be as in (iii), and suppose that $X$ is an $S$-scheme
with an action of $G:=D_S(\Gamma)$. Then the corresponding morphism
\eqref{eq_no-nam-action} and the projection $p_G:G\times_SX\to G$
induce an automorphism of the $G$-scheme $G\times_SX$, whose
composition with the projection $p_X:G\times_SX\to X$
equals \eqref{eq_no-nam-action}. We then deduce that both
\eqref{eq_no-nam-action} and $p_X$ are smooth morphisms.
This observation, together with the above correspondance between
$\Gamma$-graded algebras and $D_S(\Gamma)$-actions, is the basis
for a general method, that allows to prove properties of graded
rings, provided they can be translated as properties of the
corresponding schemes {\em which are well behaved under smooth
base changes}. We shall present hereafter a few applications of
this method.
\end{remark}

\begin{proposition}\label{prop_notnilpo}
Let $\Gamma$ be an integral monoid, $B$ a $\Gamma$-graded
(commutative, unital) ring, and suppose that the order
of any torsion element of\/ $\Gamma^\gp$ is invertible in $B$.
Then :
\begin{enumerate}
\item
$\nil(B[\Gamma])=\nil(B)\cdot B[\Gamma]$.
\item
$\nil(B)$ is a $\Gamma$-graded ideal of $B$.
\end{enumerate}
\end{proposition}
\begin{proof}(i): Clearly
$\nil(B)\cdot B[\Gamma]\subset\nil(B[\Gamma])$.
To show the converse inclusion, it suffices to prove that
$\nil(B[\Gamma])\subset\fp B[\Gamma]$ for every
prime ideal $\fp\subset B$, or equivalently that $B/\fp[\Gamma]$
is a reduced ring for every such $\fp$. Since the natural map
$\Gamma\to\Gamma^\gp$ is injective, we may further replace
$\Gamma$ by $\Gamma^\gp$, and assume that $\Gamma$ is an
abelian group. In this case, $B/\fp[\Gamma]$ is the filtered
union of the rings $B/\fp[H]$, where $H$ runs over the finitely
generated subgroups of $\Gamma$; it suffices therefore to prove
that each $B/\fp[H]$ is reduced, so we may assume that $\Gamma$
is finitely generated, and the order of its torsion subgroup is
invertible in $B$. In this case, $B/\fp[\Gamma]$ is a smooth
$B/\fp$-algebra (remark \ref{rem_etale-and-smooth}(iii)),
and the assertion follows from \cite[Ch.IV, Prop.17.5.7]{EGA4}.

(ii): The assertion to prove is that
$\nil(B)=\bigoplus_{\gamma\in\Gamma}\nil(B)\cap\gr_\gamma B$.
However, let $\theta_B:B\to B[\Gamma]$ be the map defined as
in example \ref{ex_actions-are-thetas}(ii); now, (i) implies
that $\theta_B$ restricts to a map
$\nil(B)\to\nil(B)\cdot B[\Gamma]$, whence the contention.
\end{proof}

\begin{definition}\label{def_ic}
Let $A$ be a ring, $S\subset A$ the multiplicative subset
of regular elements of $A$.

(i)\ \
For every ring homomorphism $A\to B$, we denote the
integral closure of $A$ in $B$ by :
$$
\mathrm{i.c.}(A,B).
$$

(ii)\ \
The {\em total ring of fractions} of $A$ is the ring :
$$
\Frac(A):=S^{-1}A.
$$

(iii)\ \
The {\em normalization} $A$ is $A^\nu:=\mathrm{i.c.}(A,\Frac(A))$.
We say that $A$ is {\em normal} if $A=A^\nu$.

(iv)\ \
We say that a $A$ is {\em nice} if either $A=0$ or $\Frac(A)$
has Krull dimension zero.
\end{definition}

\begin{remark}\label{rem_nice-rings}
Let $A$ be a ring.

(i)\ \
Notice that the localization map $A\to\Frac(A)$ is injective,
and every regular element of $\Frac(A)$ is invertible : the
details are left to the reader.

(ii)\ \
If $A$ has Krull dimension zero, then every regular element
of $A$ is invertible : indeed, if $a\in A$ is regular and
$\fp\in\Spec\,A$, the image of $a$ in $A_\fp$ is still regular,
hence $a\notin\fp A_\fp$, because $\fp A_\fp$ is the nilradical
of $A_\fp$. Thus, $a$ doesn't lie in any prime ideal of $A$,
whence the contention.

(iii)\ \
If $A$ is reduced and has finitely many minimal prime ideals,
then $A$ is nice. Indeed, in this case $\Frac(A)$ has also
finitely many minimal prime ideals $\fp_1,\dots,\fp_k$; now,
if $x\in\Frac(A)\setminus(\fp_1\cup\cdots\fp_k)$, and $x$ is
not invertible, then $x$ is not regular, by (i). Say that
$xy=0$ for some $y\in\Frac(A)\setminus\{0\}$; it follows that
$y\in\fp_1\cap\cdots\cap\fp_k=0$, a contradiction. So every
non-invertible element of $\Frac(A)$ lies in
$\fp_1\cup\cdots\cup\fp_k$, whence the contention.

(iv)\ \
If $A$ is noetherian and has no embedded prime ideals, then
$A$ is nice. Indeed, in this case $\Frac(A)$ fulfills the
same conditions; then, take again $x\in\Frac(A)$ in the
complement of the finitely many minimal prime ideals. If $x$
is not invertible, it is not regular according to (i), hence
it lies in some associated prime ideal of $A$
(\cite[Th.6.1(ii)]{Mat}), a contradiction.

(v)\ \
Every flat ring homomorphism $f:A\to B$ maps the regular
elements of $A$ to the regular elements of $B$, hence it
extends uniquely to a flat ring homomorphism
$$
\Frac(f):\Frac(A)\to\Frac(B).
$$
Moreover, if $f$ is injective, the same holds for $\Frac(f)$,
since the localization map $B\to\Frac(B)$ is injective by
(i) : details left to the reader.

(vi)\ \
Let $f:A\to B$ be a flat ring homomorphism, such that $A$ is
nice and $\Spec\,\kappa(\fp)\otimes_AB$ is either empty or
has Krull dimension zero for every prime ideal $\fp$ of $A$
(where $\kappa(\fp)$ denotes the residue field of $\fp$).
Then $B$ is nice and $\Frac(f)$ induces an isomorphism of
$B$-algebras :
$$
f':\Frac(A)\otimes_AB\isom\Frac(B).
$$
Indeed, let $b\in B$ be any regular element; since the induced
map $g:B\to\Frac(A)\otimes_AB$ is flat, $g(b)$ is regular in
$\Frac(A)\otimes_AB$. On the other hand, by assumption 
$\Spec\,\Frac(A)$ and the fibres of $\Spec(f)$ are either empty
or have Krull dimension zero, hence $\Frac(A)\otimes_AB$ is
either $0$ or has Krull dimension zero, so $g(b)$ is invertible,
by (ii). This shows that $g$ extends to a ring homomorphism
$\Frac(B)\to\Frac(A)\otimes_AB$ which is the inverse of $f'$,
whence the contention.

(vii)\ \
If $A$ is nice, $(A_\fp)^\nu=(A^\nu)_\fp$ for every $\fp\in\Spec\,A$.
Indeed, in light of (vi) we have :
$$
(A_\fp)^\nu=\mathrm{i.c.}(A_\fp,\Frac(A_\fp))=
\mathrm{i.c.}(A_\fp,\Frac(A)_\fp)=(A^\nu)_\fp.
$$

(viii)\ \
Especially, if $A$ is a nice and reduced ring, then $A$ is
normal if and only if $A$ is {\em locally normal}, {\em i.e.}
if and only if $A_\fp$ is a reduced local normal domain for every
$\fp\in\Spec\,A$. Indeed, the condition is sufficent, by virtue
of (vii). Conversely, suppose that $A$ is nice, reduced and normal.
In view of (vii), it suffices to check that $A_\fp$ is a domain
for every $\fp\in\Spec\,A$, and since $A_\fp$ is reduced, this
amounts to showing that $A_\fp$ has a unique minimal prime ideal.
However, $A_\fp$ is nice due to (vi), hence $X:=\Spec\,\Frac(A_\fp)$
has Krull dimension zero, and especially it is a totally
disconnected topological space. Moreover, it follows from (i)
that the inclusion map $X\to\Spec\,A_\fp$ has dense image, so
it must contain every minimal prime ideal of $\Spec\,A_\fp$, and
then $X$ is precisely the set of minimal prime ideals of $A_\fp$.
Now, if $X$ contains two distinct points, we may find non-empty
open and closed subsets $U,U'\subset X$ with $U\cap U'=\emptyset$
and $U\cup U'=X$; then there is an idempotent $e\in\Frac(A_\fp)$
which vanishes exactly on $U$. Since $A_\fp$ is normal by (vii),
it follows that $e\in A_\fp$, but then $\Spec\,A_\fp$ is the
disjoint union of two non-empty open and closed subsets, which
is absurd, since it has a unique closed point.
\end{remark}

\begin{proposition}\label{prop_integr-closure-grad}
Let $(\Gamma,+,0)$ be an integral monoid, $f:A\to B$ a morphism
of\/ $\Gamma$-graded rings, and suppose that the order of
any torsion element of\/ $\Gamma^\gp$ is invertible in $A$.
We have :
\begin{enumerate}
\item
$(\mathrm{i.c.}(A,B))[\Gamma^\gp]=
\mathrm{i.c.}(A[\Gamma^\gp],B[\Gamma^\gp])$.
\item
The grading of $B$ restricts to a $\Gamma$-grading on the
subring $\mathrm{i.c.}(A,B)$.
\item
Suppose moreover that $B$ is reduced. Then for every
$\gamma\in\Gamma$ such that $\gr_\gamma\mathrm{i.c.}(A,B)\neq 0$
there exists an integer $m>0$ with $\gr_{m\gamma}A\neq 0$.
\end{enumerate}
\end{proposition}
\begin{proof}(i): We easily reduce to the case where $\Gamma$
is finitely generated, in which case $A[\Gamma^\gp]$ is a
smooth $A$-algebra (remark \ref{rem_etale-and-smooth}(iii)),
and the assertion follows from \cite[Ch.IV, Prop.6.14.4]{EGAIV-2}.

(ii): Define $\theta_A$ and $\theta_B$ as in example
\ref{ex_actions-are-thetas}(ii), and consider the
commutative diagram :
$$
\xymatrix{
A \ar[r]^-{\theta_A} \ar[d]_f & A[\Gamma] \ar[d]^{f[\Gamma]} \\
B \ar[r]^-{\theta_B} & B[\Gamma].
}$$
Say that $b\in\mathrm{i.c.}(A,B)$;
then $\theta_B(b)\in\mathrm{i.c.}(A[\Gamma],B[\Gamma])$. In light
of (i), we deduce that
$\theta_B(b)\in(\mathrm{i.c.}(A,B))[\Gamma^\gp]\cap B[\Gamma]=
(\mathrm{i.c.}(A,B))[\Gamma]$. The claim follows easily.

(iii): Let $x$ be any non-zero element of
$\gr_\gamma\mathrm{i.c.}(A,B)$. By definition there exist
an integer $n>0$ and $a_1,\dots,a_n\in A$ such that
$x^n+a_1x^{n-1}+\cdots+a_n=0$ in $B$. For $i=1,\dots,n$
let $a'_i$ be the image of $a_i$ under the projection
$A\to\gr_{i\gamma}A$; it is easily seen that
$x^n+a'_1x^{n-1}+\cdots+a'_n=0$. If $a'_i=0$ for every
$i=1,\dots,n$, we deduce that $x^n=0$, which is absurd,
since $B$ is reduced. Thus, we have $a'_m\neq 0$ for at
least one index $m\leq n$, so $\gr_{m\gamma}A\neq 0$.
\end{proof}

\begin{definition} Let $(\Gamma,+,0)$ be an integral monoid,
$\underline A:=(A,\gr_\bullet A)$ a $\Gamma$-graded $\Z$-algebra.
Set $\gr^*_\gamma A:=\gr_\gamma A\setminus\{0\}$ for every
$\gamma\in\Gamma$, and $S^*:=\bigcup_{\gamma\in\Gamma}\gr^*_\gamma A$.

(i)\ \
We say that $\underline A$ is a {\em $*$-domain} if we have
$\gr^*_\gamma A\cdot\gr_{\gamma'}^*A\subset\gr_{\gamma+\gamma'}^*A$
for every $\gamma,\gamma'\in\Gamma$.

(ii)\ \
We say that $\underline A$ is a {\em $*$-field} if $A\neq 0$ and
every element of $S^*$ is invertible in $A$.

(iii)\ \
Notice that if $\underline A$ is a $*$-domain, then every element
of $S^*$ is regular in $A$, and obviously $S^*$ is a multiplicative
subset of $A$, so we have a unique $\Gamma^\gp$-grading on
$$
\Frac^*(\underline A):=S^{*-1}A
$$
such that the localization map $A\to\Frac^*(\underline A)$ is a
morphism of $\Gamma^\gp$-graded $\Z$-algebras.
\end{definition}

\begin{lemma}\label{lem_modules-on-*-fields}
Let $\underline A$ be a $\Gamma$-graded $*$-field, and
$\underline M:=(M,\gr_\bullet M)$ a $\Gamma$-graded
$\underline A$-module. Then $M$ is a free $A$-module,
and it admits a basis consisting of homogeneous elements
of $M$.
\end{lemma}
\begin{proof} Indeed, notice first that $\gr_0A$ is a field,
and $\Gamma_0:=\{\gamma\in\Gamma~|~\gr^*A\neq\emptyset\}$
is an abelian group (with the addition law of $\Gamma$).
Now, let $\Lambda\subset\Gamma$ be a system of representatives
for $\Gamma/\Gamma^0$ (say, with $0\in\Lambda$), and for
every $\lambda\in\Lambda$, pick a basis
$(m_{\lambda,i}~|~i\in I_\lambda)$ for the $\gr_0A$-vector space
$\gr_\lambda M$. Let us show that
$\Sigma:=(m_{\lambda,i}~|~\lambda\in\Lambda,i\in I_\lambda)$
is a basis for the $A$-module $M$. To see that $\Sigma$
is a generating family for $M$, consider any $\gamma\in\Gamma$
and any $m\in\gr_\gamma M$; then there exist $\lambda\in\Lambda$
and $\gamma'\in\Gamma_0$ such that $\gamma=\lambda+\gamma'$,
and for any $a\in\gr_{\gamma'}A$ we have a $\gr_0A$-linear
combination $a^{-1}m=\sum_{i\in J}a_im_{\lambda,i}$ for a finite
subset $J\subset I_{\lambda}$. Thus $m=\sum_{i\in J}aa_im_{\lambda,i}$,
whence the contention. Next, say that we have an $A$-linear
relation of the form :
$\sum_{\lambda\in\Lambda}\sum_{i\in I_\lambda}a_{\lambda,i}m_{\lambda,i}=0$
(where $a_{\lambda,i}=0$ for all but finitely many indices
$(\lambda,i)$); then it is easily seen that we have already
$\sum_{i\in I_\lambda}a_{\lambda,i}m_{\lambda,i}=0$ for every
$\lambda\in\Lambda$. Then it also follows straightforwardly
that $\sum_{i\in I_\lambda}\gr_\gamma(a_{\lambda,i})m_{\lambda,i}=0$
for every $\gamma\in\Gamma_0$, where $\gr_\gamma(a_{\lambda,i})$
denotes the image of $a_{\lambda,i}$ under the projection
$A\mapsto\gr_\gamma A$. For every such $\gamma$, pick
$b_\gamma\in\gr_\gamma^*A$; then we get the $\gr_0A$-linear
relation
$\sum_{i\in I_\lambda}b^{-1}_\gamma\gr_\gamma(a_{\lambda,i})m_{\lambda,i}=0$
in $\gr_\lambda M$, and therefore $\gr_\gamma(a_{\lambda,i})=0$ for
every $\lambda\in\Lambda$, every $\gamma\in\Gamma_0$ and every
$i\in I_\lambda$, {\em i.e.} $a_{\lambda,i}=0$ for all such $\lambda$
and $i$, as required.
\end{proof}

\begin{remark}\label{rem_*-domains}
(i)\ \
For every $\Gamma$-graded $*$-domain $\underline A$, the
$\Gamma^\gp$-graded $\Z$-algebra $\Frac^*(\underline A)$ is
a $*$-field, and the natural maps
$A\to\Frac^*(\underline A)\to\Frac(A)$
are injective (details left to the reader).

(ii)\ \
Let $\underline A$ and $\underline B$ be two $\Gamma$-graded
$*$-domains, and $f:\underline A\to\underline B$ an injective
morphism of $\Gamma$-graded $\Z$-algebras. Then $f$ extends
to a unique morphism of $\Gamma^\gp$-graded $\Z$-algebras
$\Frac(f)^*:\Frac^*(\underline A)\to\Frac^*(\underline B)$,
and taking into account (i) and lemma
\ref{lem_modules-on-*-fields}, we see that $\Frac^*(f)$ is
a flat and injective ring homomorphism, hence it extends in
turn uniquely to a flat and injective homomorphism of total
rings of fractions $\Frac(f):\Frac(A)\to\Frac(B)$, by remark
\ref{rem_nice-rings}(v).

(iii)\ \
Especially, let $(\underline A{}_\lambda:=
(A_\lambda,\gr^\lambda_\bullet A)~|~\lambda\in\Lambda)$ be a
filtered system of $\Gamma$-graded $\Z$-algebras, with
injective transition maps, such that $\underline A{}_\lambda$
is a $*$-domain for every $\lambda\in\Lambda$. It follows
easily that the colimit $A$ of the filtered system of
underlying rings $(A_\lambda~|~\lambda\in\Lambda)$ carries
a unique $\Gamma$-grading $\gr_\bullet$ such that the natural
ring homomorphism $j_\lambda:A_\lambda\to A$ is a morphism of
$\Gamma$-graded $\Z$-algebras
$\underline A{}_\lambda\to\underline A:=(A,\gr_\bullet)$; moreover,
$\underline A$ is also a $*$-domain, and the resulting cocone
$(\underline A{}_\lambda\to\underline A~|~\lambda\in\Lambda)$
is universal in the category of $\Gamma$-graded $\Z$-algebras
(details left to the reader). Furthermore, according to (ii),
for every $\lambda,\mu\in\Lambda$ with $\mu\geq\lambda$, the
transition map $A_\lambda\to A_\mu$ extends uniquely to a flat
and injective ring homomorphism $\Frac(A_\lambda)\to\Frac(A_\mu)$;
similarly we get a unique flat and injective extension
$\Frac(j_\lambda):\Frac(A_\lambda)\to\Frac(A)$ of $j_\lambda$,
for every $\lambda\in\Lambda$. Then it is easily seen that
the resulting cocone $(\Frac(j_\lambda)~|~\lambda\in\Lambda)$
is universal.

(iv)\ \
In the situation of (iii), suppose moreover that $A_\lambda$
is nice for every $\lambda\in\Lambda$. Then the same holds
for $A$. Indeed, from (iii) we deduce a natural identification
of $\Spec\,\Frac(A)$ with the inverse limit of the cofiltered system
of topological spaces $(\Spec\,\Frac(A_\lambda)~|~\lambda\in\Lambda)$.
But it is easily seen that the limit of a cofiltered system
of spectral spaces of dimension zero has dimension zero
(details left to the reader), whence the contention.
\end{remark}

\begin{proposition}\label{prop_*-fields}
Let $(\Gamma,+,0)$ be an abelian group,
$\underline A:=(A,\gr_\bullet A)$ a $\Gamma$-graded $*$-field,
$p\geq 1$ the characteristic exponent of the field $\gr_0A$,
and suppose that $\Gamma\setminus\{0\}$ has no $p$-torsion
elements. Then $A$ is a nice reduced normal ring.
\end{proposition}
\begin{proof} Clearly $\underline A$ is the filtered
colimit of its $\Gamma$-graded $\Z$-subalgebras
$\Delta\times_\Gamma\underline A$, where $\Delta$
ranges over the finitely generated subgroups of $\Gamma$.
Moreover, $A^\nu$ is the filtered union of its subrings
$(\Delta\times_\Gamma\underline A)^\nu$. Taking into account
remark \ref{rem_*-domains}(iv), we are thus easily reduced
to the case where $\Gamma$ is finitely generated, say
$\Gamma=\Delta_1\oplus\cdots\oplus\Delta_k$, for cyclic
abelian groups $\Delta_1,\dots,\Delta_k$. Furthermore,
after replacing $\Gamma$ by a subgroup, we may assume
that $\gr^*_\gamma\neq\emptyset$ for every $\gamma\in\Gamma$,
in which case it is easily seen that $\gr_\gamma A$ is a
one-dimensional $\gr_0A$-vector space, for every
$\gamma\in\Gamma$ (details left to the reader), and
we deduce an isomorphism of $\gr_0A$-algebras :
$$
(\Delta_1\times_\Gamma\underline A)\otimes_{\gr_0A}\cdots
\otimes_{\gr_0A}(\Delta_k\times_\Gamma\underline A)\isom A.
$$
In view of remark \ref{rem_nice-rings}(iii) we are then
reduced to checking :

\begin{claim}
If $\Gamma$ is a cyclic abelian group, $A$ is a smooth
$\gr_0A$-algebra.
\end{claim}
\begin{pfclaim}[] Let $\gamma$ be a generator of $\Gamma$;
recall that $\gr^*_\gamma A\neq\emptyset$, and pick any
$a\in\gr^*_\gamma A$. If $\Gamma\isom\Z$, then it is easily
seen that the map : $P(T)\mapsto P(a)$ yields an isomorphism
of $\gr_0A$-algebras $\gr_0A[T,T^{-1}]\isom A$. Lastly, say that
$\Gamma=\Z/n\Z$ for some $n\in\N\setminus\{0\}$; by
assumption $n$ is invertible in $\gr_0A$, and $b:=a^n\in\gr^*_0A$.
Then the same rule yields an isomorphism of $\gr_0A$-algebras
$\gr_0A[T]/(T^n-b)\isom A$.
\end{pfclaim}
\end{proof}

\begin{corollary} Let $\Gamma$ be a monoid,
$\underline A:=(A,\gr_\bullet A)$ a $\Gamma$-graded $*$-domain
and suppose that $\Gamma^\gp$ is a torsion-free abelian group.
Then $A$ is a domain.
\end{corollary}
\begin{proof} In view of remark \ref{rem_*-domains}(i) we may
replace $\underline A$ by the $\Gamma^\gp$-graded $\Z$-algebra
$\Frac^*(A)$, and assume from start that $\Gamma$ is a
torsion-free abelian group, and $\underline A$ is a $*$-field.
Then $\underline A$ is the filtered colimit of the system of
$\Gamma$-graded $\Z$-subalgebras
$\Delta\times_\Gamma\underline A$, where $\Delta$ ranges
over the finitely generated subgroups of $\Gamma$, and
it suffices to prove the assertion for such subalgebras.
Thus, we may assume $\Gamma=\Z^{\oplus r}$ for some $r\in\N$,
and after replacing $\Gamma$ by a subgroup, we may even
assume that $\gr^*_\gamma A\neq\emptyset$ for every
$\gamma\in\Gamma$. In this case, the proof of proposition
\ref{prop_*-fields} shows that $A$ is isomorphic to
$\gr_0A[T_1^{\pm 1},\dots,T_r^{\pm 1}]$, whence the contention.
\end{proof}

\begin{corollary}
Let $\Gamma$ be an integral monoid,
$\underline A:=(A,\gr_\bullet A)$ a $\Gamma$-graded $*$-domain,
$\Gamma^\nu\subset\Gamma^\gp$ the saturation of\/ $\Gamma$, and
suppose that the order of any torsion element of\/ $\Gamma^\gp$
is invertible in $A$. Then there exists a unique $\Gamma^\nu$-grading
$\gr_\bullet A^\nu$ on $A^\nu$ such that $(A^\nu,\gr_\bullet A^\nu)$
is a $\Gamma^\nu$-graded nice reduced $*$-domain, and the
inclusion map $A\to A^\nu$ is a morphism of\/
$\Gamma^\nu$-graded rings.
\end{corollary}
\begin{proof} From proposition \ref{prop_*-fields} and
remark \ref{rem_*-domains}(i) we see that
$A^\nu=\mathrm{i.c.}(A,\Frac^*(A))$. The corollary then
follows straightforwardly from proposition
\ref{prop_integr-closure-grad}.
\end{proof}

\begin{proposition}\label{prop_four-year-later}
Suppose that $(\Gamma',+)\to(\Gamma,+)$ is a morphism
of fine monoids, $B$ a finitely generated (resp. finitely
presented) $\Gamma$-graded $R$-algebra, and $M$ a finitely
generated (resp. finitely presented) $\Gamma$-graded $B$-module.
We have :
\begin{enumerate}
\item
$\Gamma'\times_\Gamma B$ is a finitely generated (resp.
finitely presented) $R$-algebra.
\item
$\Gamma'\times_\Gamma M$ is a finitely generated (resp. finitely
presented) $\Gamma'\times_\Gamma B$-module.
\item
$\gr_\gamma M$ is a finitely generated (resp. finitely presented)
$\gr_0B$-module, for every $\gamma\in\Gamma$.
\end{enumerate}
\end{proposition}
\begin{proof} Let $B_M$ denote the direct sum $B\oplus M$,
endowed with the $R$-algebra structure given by the rule :
$$
(b_1,m_1)\cdot(b_2,m_2)=(b_1b_2,b_1m_2+b_2m_1)
\qquad
\text{for every $b_1,b_2\in B$ and $m_1,m_2\in M$}.
$$
Notice that $B_M$ is characterized as the unique $R$-algebra
structure for which $M$ is an ideal with $M^2=0$, the natural
projection $B_M\to B$ is a map of $R$-algebras, and the
$B$-module structure on $M$ induced via $\pi$ agrees with
the given $B$-module structure on $M$.

\begin{claim}\label{cl_B_M-is-here}
The following conditions are equivalent :
\begin{enumerate}
\alphaenu
\item
The $R$-algebra $B_M$ is finitely generated (resp. finitely presented).
\item
$B$ is a finitely generated (resp. finitely presented) $R$-algebra
and $M$ is a finitely generated (resp. finitely presented) $B$-module.
\end{enumerate}
\end{claim}
\begin{pfclaim}(b)$\Rightarrow$(a): Suppose first that $B$ is
a finitely generated $R$-algebra, and $M$ is a finitely generated
$B$-module. Pick a system of generators $\Sigma_B:=\{b_1,\dots,b_s\}$
for $B$ and $\Sigma_M:=\{m_1,\dots,m_k\}$ for $M$. Then it is easily
seen that $\Sigma_B\cup\Sigma_M$ generates the $R$-algebra $B_M$.

For the finitely presented case, pick a surjection of $R$-algebras
$\phi:R[T_1,\dots,T_s]\to B$ and of $B$-module $\psi:B^{\oplus k}\to M$.
Let $\Sigma'_B$ be a finite system of generators of the ideal
$\Ker\,\phi$. Pick also a finite system $b_1,\dots,b_r$ of generators
of the $B$-module $\Ker\,\psi$; we may write
$b_i=\sum_{j=1}^kb_{ij}e_j$ for certain $b_{ij}\in B$ (where
$e_1,\dots,e_k$ is the standard basis of $B^{\oplus k}$). For
every $i\leq r$ and $j\leq k$, pick $P_{ij}\in\phi^{-1}(b_{ij})$.
It is easily seen that $B_M$ is isomorphic to the $R$-algebra
$R[T_1,\dots,T_{s+k}]/I$, where $I$ is generated by 
$\Sigma'_B\cup\{\sum_{j=1}^kP_{ij}T_{j+s}~|~i=1,\dots,r\}
\cup\{T_{i+s}T_{j+s}~|~0\leq i,j\leq k\}$.

(a)$\Rightarrow$(b): Suppose that $B_M$ is a finitely generated
$R$-algebra, and let $c_1,\dots,c_n$ be a system of generators.
For every $i=1,\dots,n$, we may write $c_i=b_i+m_i$ for unique
$b_i\in B$ and $m_i\in M$. Since $M^2=0$, it is easily seen that
$m_1,\dots,m_n$ is a system of generators for the $B$-module $M$,
and clearly $b_1,\dots,b_n$ is a system of generators for the
$R$-algebra $B$.

Next, suppose that $B_M$ is finitely presented over $R$.
We may find a system of generators
$b_1,\dots,b_s,m_1,\dots,m_k$ of $B_M$ with $b_i\in B$ and
$m_j\in M$ for every $i\leq s$ and $j\leq k$. We deduce a
surjection of $R$-algebras
$$
\phi:R[T_1,\dots,T_{s+k}]/(T_{s+i}T_{s+j}~|~0\leq i,j\leq k)
\to B_M
$$
such that $T_i\mapsto b_i$ for every $i\leq s$ and
$T_{j+s}\mapsto m_j$ for every $j\leq k$. It is easily
seen that $\Ker\,\phi$ is generated by the classes of
finitely many polynomials $P_1,\dots,P_r$, where
$$
P_i=Q_i(T_1,\dots,T_s)+\sum_{j=1}^kT_{s+j}Q_{ij}(T_1,\dots,T_s)
\qquad
i=1,\dots,r
$$
for certain polynomials $Q_i,Q_{ij}\in R[T_1,\dots,T_s]$.
It follows easily that $B=R[T_1,\dots,T_s]/I$, where $I$
is the ideal generated by $Q_1,\dots,Q_r$, and $M$ is
isomorphic to the $B$-module $B^{\oplus k}/N$, where
$N$ is the submodule generated by the system
$\{\sum_{j=1}^kQ_{ij}(b_1,\dots,b_s)e_j~|~i=1,\dots,r\}$
\end{pfclaim}

Suppose now that $B$ is a finitely generated $R$-algebra;
then $B$ is generated by finitely many homogeneous elements, say
$b_1,\dots,b_s$ of degrees respectively $\gamma_1,\dots,\gamma_s$.
Thus, we may define surjections of monoids
\set\begin{equation}\label{eq_viavia}
\N^{\oplus s}\to\Gamma
\quad :\quad
e_i\mapsto\gamma_i
\qquad
\text{for $i=1,\dots,s$}
\end{equation}
(where $e_1,\dots,e_s$ is the standard basis of $\N^{\oplus s}$)
and of $R$-algebras $\phi:C\to B$, where $C:=R[\N^{\oplus s}]$ is a
free polynomial $R$-algebra. Notice that $C$ is a
$\N^{\oplus s}$-graded $R$-algebra, and via \eqref{eq_viavia}
we may regard $\phi$ as a morphism of $\Gamma$-graded $R$-algebras
$C_{/\Gamma}\to B$. Then $I:=\ker\,\phi$ is a $\Gamma$-graded
ideal of $C$, and if we set $P:=\N^{\oplus s}\times_\Gamma\Gamma'$
we deduce an isomorphism of $\Gamma'$-graded $R$-algebras
$$
B':=\Gamma'\times_\Gamma B\isom
(P\times_{\N^{\oplus s}}C)_{/\Gamma'}/(\Gamma'\times_\Gamma I)
$$
(see remark \ref{rem_etale-and-smooth}(i)).

\begin{claim}\label{cl_cucu-again}
$P\times_{\N^{\oplus s}}C$ is a finitely presented $R$-algebra.
\end{claim}
\begin{pfclaim} Indeed, this $R$-algebra is none else than
$R[P]$, hence the assertion follows from corollary
\ref{cor_fibres-are-fg} and lemma \ref{lem_finite-pres}(i).
\end{pfclaim}

From claim \ref{cl_cucu-again} it follows already that $B'$ is a
finitely generated $R$-algebra. Now, suppose that $M$ is
a finitely generated $B$-module, and set
$M':=\Gamma'\times_\Gamma M$; notice that
\set\begin{equation}\label{eq_shurik}
\Gamma'\times_\Gamma(B_M)=B'_{M'}.
\end{equation}
In view of claim \ref{cl_B_M-is-here}, we deduce that $M'$
is a finitely generated $B'$-module. Next, in case $B$ is a
finitely presented $R$-algebra, $I$ is a finitely generated
ideal of $C_{/\Gamma}$; as we have just seen, this implies
that $\Gamma'\times_\Gamma I$ is a finitely generated
$(\Gamma'\times_\Gamma C_{/\Gamma})$-module, and then claim
\ref{cl_cucu-again} shows that $B'$ is a finitely presented
$R$-algebra. This concludes the proof of (i).

Lastly, if moreover $M$ is a finitely presented
$B$-module, assertion (i), together with \eqref{eq_shurik}
and claim \ref{cl_B_M-is-here} say that $M'$ is a finitely
presented $B'$-module; thus, also assertion (ii) is proven.

(iii): For any given $\gamma\in\Gamma$, let us consider the
morphism $f:\N\to\Gamma$ such that $1\mapsto\gamma$, and set
$B':=\N\times_\Gamma B$. By (i), the $R$-algebra $B'$
is finitely generated (resp. finitely presented), and the
$B'$-module $M':=\N\times_\Gamma M$ is finitely generated
(resp. finitely presented). After replacing $B$ by $B'$
and $M$ by $M'$, we may then assume from start that $\Gamma=\N$,
and we are reduced to showing that $\gr_1M$ is a finitely generated
(resp. finitely presented) $\gr_0B$-module.

Let $m_1,\dots,m_t$ be a system of generators of $M$ consisting
of homogeneous elements of degrees respectively $j_1,\dots,j_t$.
We endow $B^{\oplus t}$ with the $\N$-grading such that
$$
\gr_k B^{\oplus t}:=\bigoplus_{i=1}^t \gr_{k-j_i}Be_i
$$
(where $e_1,\dots,e_t$ is the standard basis of $B^{\oplus t}$);
then the $B$-linear map $B^{\oplus t}\to M$ given by the rule
$e_i\mapsto m_i$ for every $i=1,\dots,t$ is a morphism of
$\N$-graded $B$-modules, and if $M$ is finitely presented,
its kernel is generated by finitely many homogeneous elements
$b_1,\dots,b_s$. In the latter case, endow again $B^{\oplus s}$
with the unique $\N$-grading such that the $B$-linear map
$\phi:B^{\oplus s}\to B^{\oplus t}$ given by the rule
$e_i\mapsto b_i$ for every $i=1,\dots,s$ is a morphism of
$\N$-graded $B$-modules. Now, in order to check that $\gr_1M$
is a finitely generated $\gr_0B$-module, it suffices to show
that the same holds for $\gr_1B^{\oplus t}$. The latter is
a direct sum of $\gr_0B$-modules isomorphic to either $\gr_0B$
or $\gr_1B$. Likewise, if $M$ is finitely presented,
$\gr_1M=\Coker\,\gr_1\phi$, and again, $\gr_1B^{\oplus s}$ is
a direct sum of $\gr_0B$-modules isomorphic to either $\gr_0B$ or
$\gr_1B$; hence in order to check that $\gr_1M$ is a finitely
presented $\gr_0B$-module, it suffices to show that $\gr_1B$ is
a finitely presented $\gr_0B$-module.  In either event, we are
reduced to the case where $\Gamma=\N$ and $M=B$.

However, from (ii) we deduce especially that
$\gr_0B=\{0\}\times_\N B$ is a finitely generated (resp. finitely
presented) $R$-algebra, hence $B$ is a finitely generated
(resp. finitely presented) $B_0$-algebra as well; we may then
assume that $R=\gr_0B$. Let $\Sigma$ be a system of homogeneous
generators for the $R$-algebra $B$; we may then also assume that
\set\begin{equation}\label{eq_no-intersect}
\Sigma\cap\gr_0B=\emptyset.
\end{equation}
Then it is easily seen that the $R$-module $\gr_1B$ is generated
by $\Sigma\cap\gr_1B$. Lastly, if $B$ is a finitely presented
$R$-algebra, we consider the natural surjection
$\psi:R[\Sigma]\to B$ from the free polynomial
$R$-algebra generated by the set $\Sigma$, and endow $R[\Sigma]$
with the unique grading for which $\psi$ is a map of $\N$-graded
$R$-algebras; then $I:=\Ker\,\psi$ is a finitely generated
$\N$-graded ideal with $\gr_0I=0$. As usual, we pick a finite
system $\Sigma'$ of homogeneous generators for $I$; clearly
$B_1$ is isomorphic to $\gr_1 B_0[\Sigma]/\gr_1I$. On the other
hand, \eqref{eq_no-intersect} easily implies that $\gr_1R[\Sigma]$
is a free $R$-module of finite rank, and moreover $\gr_1I$ is
generated by $\Sigma'\cap\gr_1I$; especially, $\gr_1B$ is a
finitely presented $R$-module in this case, and the proof is
complete.
\end{proof}

\sset\subsubsection{}\label{subsec_shift-by-Gamma}
Let $(\Gamma,+)$ be a monoid, $R$ a ring,
$\underline B:=(B,\gr_\bullet B)$ a $\Gamma$-graded
$R$-algebra, and $M$ a $\Gamma$-graded $\underline B$-module.
We denote by $M[\gamma]$ the $\Gamma$-graded
$\underline B$-module whose underlying $B$-module is $M$,
and whose grading is given by the rule :
$$
\gr_\beta M[-\gamma]:=\bigoplus_{\delta+\gamma=\beta}\gr_\delta M
\qquad
\text{for every $\gamma\in\Gamma$}.
$$

\begin{remark}\label{rem_graded-resol}
(i)\ \ 
In the situation of \eqref{subsec_shift-by-Gamma}, pick any
system $\bx:=(x_i~|~i\in I)$ of homogeneous generators of $M$, and
say that $x_i\in\gr_{\gamma_i}M$ for every $i\in I$. Then we may
define a surjective map of $\Gamma$-graded $\underline B$-modules
$$
L:=\bigoplus_{i\in I}B[-\gamma_i]\to M
\quad : \quad
e_i\mapsto x_i
\qquad
\text{for every $i\in I$}
$$
where $(e_i~|~i\in I)$ denotes the canonical basis of
the free $B$-module $L$ (notice that $e_i\in\gr_{\gamma_i}L$
for every $i\in I$).

(ii)\ \
In case $M$ is a finitely generated $B$-module, we may
pick a finite system $\bx$ as above, and then $L$ shall
be a free $B$-module of finite rank.

(iii)\ \
Especially, suppose that $B$ is a coherent ring and $M$ is
finitely presented as a $B$-module; then, in the situation
of (ii), the kernel of the surjection $L\to M$ shall be again
a finitely presented $\Gamma$-graded $\underline B$-module,
so we can repeat the above construction, and find inductively
a resolution
$$
\Sigma
\quad :\quad
\cdots\to L_n\xrightarrow{\ d_n\ }L_{n-1}\to\cdots\to L_0
\xrightarrow{\ d_0\ }M
$$
such that $L_n$ is a free $B$-module of finite rank,
and the map $d_n$ is a morphism of $\Gamma$-graded
$\underline B$-modules, for every $n\in\N$.

(iv)\ \
In the situation of (iii), suppose furthermore that
$B$ is a flat $R$-algebra, in which case $\gr_\gamma B$
is a flat $R$-module, for every $\gamma\in\Gamma$. Then
it is clear that the resolution $\Sigma$ yields, in each
degree $\gamma\in\Gamma$ a flat resolution $\Sigma_\gamma$
of the $R$-module $\gr_\gamma M$.
\end{remark}

\subsection{Differential graded algebras}
\label{subsec_diff-grad-alg}
The material of this paragraph shall be applied in section
\ref{sec_homotopy}, in order to study certain strictly
anti-commutative graded algebras constructed via homotopical
algebra. Especially, the graded algebras appearing in
this paragraph are usually {\em not} commutative, unless
it is explicitly said otherwise.

\begin{definition}\label{def_anti-commutative}
(i)\ \
Let $A$ be any ring, and $\underline B:=(B,\gr_\bullet B)$
a $\Z$-graded (associative, unital) $A$-algebra.
We say that $\underline B$ is {\em strictly anti-commutative}
(or {\em alternating}), if we have
\begin{enumerate}
\alphaenu
\item
$a\cdot b=(-1)^{pq}\cdot b\cdot a$ for every $p,q\in\Z$,
every $a\in\gr_pB$, and every $b\in\gr_qB$
\item
$a\cdot a=0$ for every $p\in\Z$ and every $a\in\gr_{2p+1}B$.
\end{enumerate}
If only condition (a) holds, we say that $\underline B$ is
{\em anti-commutative}.

(ii)\ \
Let $\underline M:=(M,\gr_\bullet M)$ be any $\Z$-graded
$\underline B$-bimodule. An {\em $A$-linear graded
derivation} from $\underline B$ to $\underline M$ is a morphism
of graded $A$-modules $\partial:\underline B\to\underline M$
such that
$$
\partial(xy)=\partial(x)\cdot y+x\cdot\partial(y)
\qquad
\text{for every $x,y\in B$}.
$$
\end{definition}

\begin{remark}\label{rem_alternating-alg}
(i)\ \
Let $A\tdu\Alt$ be the full subcategory of the category
of $\Z$-graded $A$-algebras, whose objects are the
strictly anti-commutative graded $A$-algebras. Define
also the category $A\tdu\AlgMod$ as in \cite[Def.2.5.22]{Ga-Ra}.
We have an obvious functor
$$
\gr_{0,1}:A\tdu\Alt\to A\tdu\AlgMod
\qquad
(B,\gr_\bullet B)\mapsto(\gr_0B,\gr_1B)
$$
and it is easily seen that the functor $\gr_{0,1}$ is right
adjoint to the functor
$$
\Lambda^\bullet_\bullet:A\tdu\AlgMod\to A\tdu\Alt
\qquad
(B,M)\mapsto\Lambda^\bullet_BM
$$
that assigns to every object $(B,M)$ the graded $A$-algebra
of the $B$-linear exterior powers of $M$; the details shall
be left to the reader.

(ii)\ \
For every pair of graded $A$-algebras
$(B,\gr_\bullet B),(B',\gr_\bullet B')$ we define
a {\em tensor product}
$$
(B'',\gr_\bullet B''):=(B,\gr_\bullet B)\otimes_A(B',\gr_\bullet B')
$$
which is the graded $A$-algebra with $B'':=B\otimes_AB'$,
and $\gr_nB'':=\bigoplus_{i+j=n}\gr_iB\otimes_A\gr_jB'$
for every $n\in\Z$. The multiplication law of $B''$ is the
direct sum of the maps
$$
(\gr_iB\otimes_A\gr_jB')\times(\gr_kB\otimes_A\gr_lB')
\to\gr_{i+k}B\otimes_A\gr_{j+l}B'
\qquad
(b\otimes b',c\otimes c')\mapsto(-1)^{jk}bc\otimes b'c'.
$$
It is easily seen that if $(B,\gr_\bullet B)$ and
$(B',\gr_\bullet B')$ are anti-commutative (resp. alternating),
the same holds for their tensor product. Moreover, if
$(B,\gr_\bullet B)$ and $(B',\gr_\bullet B')$ are alternating,
then $(B'',\gr_\bullet B'')$ represents the coproduct of
$(B,\gr_\bullet B)$ and $(B',\gr_\bullet B')$ in the category
$A\tdu\Alt$.
Again, we leave the details to the reader.

(iii)\ \
Likewise, all finite coproducts are representable in
$A\tdu\AlgMod$ : namely, the coproduct of any two objects
$(B,M),(B',M')$ is represented by :
$$
(B'',M''):=(B\otimes_AB',(M\otimes_AB')\oplus(B\otimes_AM'))
$$
with the obvious universal cocone
$(B,M)\to(B'',M'')\leftarrow(B',M')$. By (i), (ii), and
proposition \ref{prop_was-get-maddd}(iv), we deduce a natural
isomorphism of alternating graded $A$-algebras :
$$
\Lambda_{B''}^\bullet M''\isom
\Lambda^\bullet_BM\otimes_A\Lambda^\bullet_{B'}M'.
$$
\end{remark}

\begin{definition}\label{def_diff-graded-alg}
Let $A$ be any ring.

(i)\ \
A {\em differential graded $A$-algebra} is the datum of
\begin{itemize}
\item
a complex $(B^\bullet,d_B^\bullet)$ of $A$-modules
\item
an $A$-linear map $\mu^{pq}:B^p\otimes_AB^q\to B^{p+q}$
for every $p,q\in\Z$
\end{itemize}
such that the following holds :
\begin{enumerate}
\alphaenu
\item
Set $B:=\bigoplus_{p\in\Z}B^p$; then the system of maps
$\mu^{\bullet\bullet}$ adds up to a map $\mu:B\otimes_AB\to B$,
and one requires that the resulting pair $(B,\mu)$ is
an associative unital $\Z$-graded $A$-algebra. Then, one
sets $a\cdot b:=\mu(a\otimes b)$, for every $a,b\in B$.
\item
We have the identities
$$
d_B^{p+q}(a\cdot b)=d_B^p(a)\cdot b+(-1)^p\cdot a\cdot d_B^q(b)
\qquad
\text{for every $p,q\in\Z$ and every $a\in B^p$, $b\in B^q$}.
$$
\end{enumerate}
We call $B$ the {\em graded $A$-algebra associated to}
$B^\bullet$. Then we also say that $B^\bullet$ is
{\em anti-commutative} (resp. {\em strictly anti-commutative})
if the same holds for $B$ (definition \ref{def_anti-commutative}).

(ii)\ \
A {\em morphism $B^\bullet\to C^\bullet$ of differential graded
$A$-algebras} is a map of complexes of $A$-modules such that
the induced map of associated graded $A$-modules is a map of
$A$-algebras. We denote the resulting category of differential
graded $A$-algebras by :
$$
A\dga.
$$
We shall also consider the full subcategory of $A\dga$ denoted
$$
A\dgalt
$$
whose objects are the strictly anti-commutative differential
graded $A$-algebras.
 
(iii)\ \
Let $(B^\bullet,d_B^\bullet)$ be a differential graded $A$-algebra,
$B$ its associated graded $A$-algebra, and $(M^\bullet,d^\bullet_M)$
a complex of $A$-modules. 
\begin{enumerate}
\alphaenu
\item
We say that $M^\bullet$ is a {\em left $B^\bullet$-module} if
the $A$-module $M:=\bigoplus_{p\in\Z}M^p$ is a graded left
$B$-module (for the natural $\Z$-grading on $M$), and we have
$$
d^{p+q}_M(b\cdot m)=(d_B^pb)\cdot m+(-1)^p\cdot b\cdot d^q_M(m)
$$
for every $p,q\in\Z$, every $b\in B^p$, and every $m\in M^q$.
\item
We say that $M^\bullet$ is a {\em right $B^\bullet$-module} if
$M$ is a graded right $B$-module, and we have
$$
d^{p+q}_M(m\cdot a)=d^q_M(m)\cdot a+(-1)^q\cdot m\cdot d^p_M(a)
$$
for every $p,q\in\Z$, every $b\in B^p$, and every $m\in M^q$.
\item
We say that $M^\bullet$ is a {\em $B$-bimodule} if it is both
a left and right $B^\bullet$-module, and with these
$B^\bullet$-modules structures, the $A$-module $M$ becomes
a $B$-bimodule ({\em i.e.} the left multiplication commutes
with the right multiplication).
\end{enumerate}
We call $M$ the {\em graded $B$-module associated to} $M^\bullet$.
A morphism $M^\bullet\to N^\bullet$ of left (resp. right, resp.
bi-) $B^\bullet$-modules is a map of complexes of $A$-modules,
such that the induced map of associated graded $A$-modules is
a map of left (resp. right, resp. bi-) $B$-modules.
\end{definition}

\begin{remark}\label{rem_dgas}
Let $B^\bullet$ be any differential graded $A$-algebra.

(i)\ \
Notice that condition (b) of definition \ref{def_diff-graded-alg}(i)
is the same as saying that $\mu^{\bullet\bullet}$ induces a map of
complexes
$$
\mu^\bullet:B^\bullet\otimes_AB^\bullet\to B^\bullet.
$$
Moreover, in light of example \ref{ex_monoidal}(i), we see
that $B^\bullet$ is anti-commutative if and only if the diagram
$$
\xymatrix{ B^\bullet\otimes_AB^\bullet \ar[rd]_{\mu^\bullet}
\ar[rr]^-\sim & & B^\bullet\otimes_AB^\bullet \ar[ld]^{\mu^\bullet} \\
& B^\bullet
}$$
commutes, where the horizontal arrow is the isomorphism
\eqref{eq_commut-constraint} that swaps the two tensor factors.
Hence, in some sense this is actually a commutativity condition.

(ii)\ \
Likewise, if $M^\bullet$ is a complex of $A$-modules with
a graded left (right) $B$-module structure on the associated
graded $A$-module $M$, then $M^\bullet$ is a left (resp. right)
$B^\bullet$-module if and only if the scalar multiplication
of the $B$-module $M$ induces a map of complexes
$$
B^\bullet\otimes_AM^\bullet\to M^\bullet
\qquad
\text{(\ resp. $M^\bullet\otimes_AB^\bullet\to M^\bullet$\ )}.
$$

(iii)\ \
It is easily seen that the multiplication maps $\mu^{pq}$
induce $A$-linear maps
$$
(H^pB^\bullet)\otimes_A(H^qB^\bullet)\to H^{p+q}B^\bullet
\qquad
\text{for every $p,q\in\Z$}
$$
anf if we let $H^\bullet B^\bullet:=\bigoplus_{p\in\Z}H^pB^\bullet$,
then the resulting map
$$
(H^\bullet B^\bullet)\otimes_A(H^\bullet B^\bullet)\to
H^\bullet B^\bullet
$$
endows $H^\bullet B^\bullet$ with a structure of $\Z$-graded
associative unital $A$-algebra, which shall be strictly
anti-commutative whenever the same holds for $B^\bullet$.
Likewise, if $M^\bullet$ is any left (resp. right, resp. bi)
$B^\bullet$-module, then $H^\bullet M^\bullet$ is naturally a
$\Z$-graded left (resp. right, resp. bi) $H^\bullet B^\bullet$-module.

(iv)\ \
Let $M^\bullet$ be a left $B^\bullet$-module, and denote
by $\mu_M^{pq}:B^p\otimes_AM^q\to M^{p+q}$ the $(p,q)$-graded
component of the scalar multiplication of $M^\bullet$, for
every $p,q\in\Z$. Then $M^\bullet[1]$ is also naturally
a left $B^\bullet$-module, with scalar multiplication
$\mu^{\bullet\bullet}_{M[1]}$ given by the rule :
$$
\mu^{pq}_{M[1]}:=(-1)^p\cdot\mu_M^{p,q+1}
\qquad\text{for every $p,q\in\Z$}.
$$
Likewise, if $N^\bullet$ is a right $B^\bullet$-module,
with scalar multiplication $\mu_N^{\bullet\bullet}$, then
$N^\bullet[1]$ is naturally a right $B^\bullet$-module, with
multiplication $\mu_{N[1]}^{\bullet\bullet}$ given by the
rule :
$$
\mu^{pq}_{N[1]}:=\mu_N^{p,q+1}
\qquad\text{for every $p,q\in\Z$}.
$$
Lastly, if $P^\bullet$ is a $B^\bullet$-bimodule, then the
left and right $B^\bullet$-module structure defined above
on $P^\bullet[1]$, amount to a natural $B$-bimodule structure
on $P^\bullet[1]$.

(v)\ \
Let $C^\bullet$ be any $\Z$-graded $A$-algebra, $N^\bullet$
any $\Z$-graded left (resp. right, resp. bi-) $B$-module.
We let $N[1]^\bullet$ be the graded $A$-module given by the
rule $N[1]^p:=N^{p+1}$ for every $p\in\Z$. We shall always
view $N[1]^\bullet$ as a left (resp. right, resp. bi-) $B$-module,
via the scalar multiplications obtained from those of
$N^\bullet$, following the rules spelled out in (iv). This
ensures that the functor from $B^\bullet$-modules to
$H^\bullet B^\bullet$-modules that assigns to any
$B^\bullet$-module $M^\bullet$ its homology
$H^\bullet M^\bullet$, is compatible with shift operators.
\end{remark}

\sset\subsubsection{}\label{subsec_finally-got-it}
Let $A$ be ring, $(B^\bullet,d_B^\bullet)$ any differential
graded $A$-algebra, and $I^\bullet\subset B^\bullet$ a
(graded) two-sided ideal of $B^\bullet$ ({\em i.e.} a
bi-submodule of the $B^\bullet$-bimodule $B^\bullet$). Let
$$
\partial:H^\bullet(B^\bullet/I^\bullet)\to H^\bullet I^\bullet[1]
$$
denote the natural map induced by the short exact sequence of
complexes
$$
0\to I^\bullet\to B^\bullet\to B^\bullet/I^\bullet\to 0.
$$
We have :

\begin{lemma}\label{lem_finally-got-it}
In the situation of \eqref{subsec_finally-got-it}, the
map $\partial$ is an $A$-linear graded derivation of the
graded $A$-algebra $H^\bullet(B^\bullet/I^\bullet)$.
\end{lemma}
\begin{proof} Indeed, let $a$ and $b$ be any two cycles
of the complex $B^\bullet/I^\bullet$ in degree $p$ and $q$,
and $\bar a$, $\bar b$ the respective classes; lift $a$
and $b$ to some elements $\tilde a\in B^p$ and
$\tilde b\in B^q$, so that $\partial(\bar a\cdot\bar b)$
is the class in $H^{p+q+1}I^\bullet$ of
$d_B(\tilde a\cdot\tilde b)=
d_B(\tilde a)\cdot\tilde b+(-1)^p\cdot\tilde a\cdot d_B(\tilde b)$.
Since $d_B(\tilde a)$ (resp. $d_B(\tilde b)$) represents
the class of $\partial(\bar a)$ (resp. of $\partial(\bar b)$),
the assertion follows from the explicit description of the
bimodule structure on $I^\bullet[1]$ provided by remark
\ref{rem_dgas}(iv).
\end{proof}

\begin{remark}\label{rem_diff-products}
Let $B^\bullet$ and $C^\bullet$ be two differential graded
$A$-algebras, and denote by $B$ and $C$ the respective
associated graded $A$-algebras. As in remark
\ref{rem_alternating-alg}(ii), we may endow $B\otimes_AC$
with a natural structure of graded $A$-algebra. In terms
of morphisms of complexes, the multiplication law of
$B\otimes_AC$ then corresponds to the composition
$$
(B^\bullet\otimes_AC^\bullet)\otimes_A(B^\bullet\otimes_AC^\bullet)
\isom(B^\bullet\otimes_AB^\bullet)\otimes_A(C^\bullet\otimes_AC^\bullet)
\xrightarrow{\ \mu^\bullet_B\otimes_A\mu^\bullet_C\ }
B^\bullet\otimes_AC^\bullet
$$
where the first isomorphism is obtained by composing the
associativity isomorphisms of example \ref{ex_monoidal}(ii)
and the isomorphisms \eqref{eq_commut-constraint} that swap
the tensor factors. Here $\mu^\bullet_B$ and $\mu^\bullet_C$
are the multiplication maps of $B^\bullet$ and $C^\bullet$.
Thus, we obtain on the complex $B^\bullet\otimes_AC^\bullet$
a natural structure of differential graded $A$-algebra,
and taking into account remark \ref{rem_alternating-alg}(ii)
it is easily seen that if $B^\bullet$ and $C^\bullet$ are strictly
anti-commutative, then $B^\bullet\otimes_AC^\bullet$ represents
the coproduct of $B^\bullet$ and $C^\bullet$ in $A\dgalt$.
\end{remark}

\sset\subsubsection{}
Suppose now that both $B^\bullet$ and $C^\bullet$ are bounded
above complexes; presumably, in this case there is a canonical
way to define a differential graded algebra structure as well
on (suitable representatives for)
$$
D^\bullet:=B^\bullet\derotimes_AC^\bullet
$$
in such a way that this structure is well defined as an object
of the derived category of $A\dga$ (the latter should be, as
usual, the localization of $A\dga$, by the multiplicative system
of quasi-isomorphisms). More modestly, we shall endow the graded
$A$-module $H_\bullet D^\bullet$ with a natural structure of
associative graded $A$-algebra, in such a way that the natural map
\set\begin{equation}\label{eq_from-D-to-BAC}
H_\bullet D^\bullet:=\bigoplus_{i\in\Z}\Tor^A_i(B^\bullet,C^\bullet)
\to H_\bullet(B^\bullet\otimes_AC^\bullet)
\end{equation}
is a morphism of graded $A$-algebras. The multiplication
of $H^\bullet D^\bullet$ is defined as the composition
$$
H_iD^\bullet\otimes_AH_jD^\bullet\xrightarrow{\ \alpha\ }
\Tor^A_{i+j}(B^\bullet\otimes_AB^\bullet,C^\bullet\otimes_AC^\bullet)
\xrightarrow{\ \mu\ } H_{i+j}D^\bullet
\qquad
\text{for every $i,j\in\Z$}
$$
where $\alpha$ is the bilinear pairing provided by
\eqref{subsec_pairing-Tors}, and with
$\mu:=\Tor^A_{i+j}(\mu^\bullet_B,\mu^\bullet_C)$.

Let us check first that the foregoing rule does define an
associative multiplication on $H_\bullet D^\bullet$. To this
aim, set $B^\bullet_2:=B^\bullet\otimes_AB^\bullet$,
$B^\bullet_3:=B^\bullet\otimes_AB^\bullet_2$ and define
likewise $C^\bullet_2$ and $C^\bullet_3$; a little diagram
chase, together with \eqref{eq_laborious}, reduces to
verifying the commutativity of the diagram
$$
\xymatrix{
\Tor_i^A(B^\bullet_2,C^\bullet_2)\otimes_AH_jD^\bullet
\ar[rr]^-{\mu\otimes_A\one_{H_iD^\bullet}} \ar[d] & &
H_iD^\bullet\otimes_AH_jD^\bullet \ar[d] & &
H_iD^\bullet\otimes_A\Tor_j^A(B^\bullet_2,C^\bullet_2)
\ar[ll]_-{\one_{H_jD^\bullet}\otimes_A\mu} \ar[d] \\
\Tor_{i+j}^A(B^\bullet_3,C^\bullet_3)
\ar[rr]^-\gamma & & \Tor_{i+j}^A(B^\bullet_2,C^\bullet_2)
& & \ar[ll]_-\delta \Tor_{i+j}^A(B^\bullet_3,C^\bullet_3)
}$$
whose vertical arrows are the bilinear pairings of
\eqref{subsec_pairing-Tors}, and with
$$
\gamma:=\Tor^A_{i+j}(\mu_B^\bullet\otimes_A\one_{B^\bullet},
\mu_C^\bullet\otimes_A\one_{C^\bullet})
\qquad
\delta:=\Tor^A_{i+j}(\one_{B^\bullet}\otimes_A\mu_B^\bullet,
\one_{C^\bullet}\otimes_A\mu_C^\bullet).
$$
We show the commutativity of the left subdiagram;
the same argument applies to the right one.
Unwinding the definitions, we come down to checking the
commutativity, in $\sD(A\Mod)$, of the diagram of
complexes
$$
\xymatrix{
(P^\bullet_{B_2}\otimes_AC_2^\bullet)\otimes_A(P^\bullet_B\otimes_AC^\bullet)
\ar[r]^-\sim
\ar[d]_{\theta^\bullet} &
(P^\bullet_{B_2}\otimes_AP^\bullet_B)\otimes_AC^\bullet_3
\ar[rrr]^-{\phi^\bullet_{12,3}\otimes_A\one_{C^\bullet_3}} & & &
P^\bullet_{B_3}\otimes_AC^\bullet_3
\ar[d]^{\eta^\bullet} \\
(P^\bullet_B\otimes_AC^\bullet)\otimes_A(P^\bullet_B\otimes_AC^\bullet)
\ar[r]^-\sim &
(P^\bullet_B\otimes_AP^\bullet_B)\otimes_AC^\bullet_2
\ar[rrr]^-{\phi^\bullet_{1,23}\otimes_A\one_{C^\bullet_3}} & & &
P^\bullet_{B_2}\otimes_AC^\bullet_2
}$$
(notation of \eqref{subsec_give_Ps}), with $\theta^\bullet:=
(P^\bullet_{\mu_B}\otimes_A\mu_C^\bullet)\otimes_A\one_{P^\bullet_B\otimes_AC^\bullet}$
and $\eta^\bullet:=
P^\bullet_{\mu_B\otimes_A\one_B}\otimes_A(\mu_C^\bullet\otimes_A\one_{C^\bullet})$,
and where $\phi^\bullet_{12,3}$ and $\phi^\bullet_{1,23}$ are as in
\eqref{subsec_pairing-Tors} (and with the isomorphisms given
by compositions of associativity and swapping isomorphisms).
A direct inspection shows that this diagram commutes up
to homotopy, as required.

Next, we check that \eqref{eq_from-D-to-BAC} is a
map of graded $A$-algebras. Unwinding the definitions, we see
that -- with the notation of \eqref{subsec_pairing-Tors} --
the multiplication of $H_\bullet D^\bullet$ is the map
on homology induced by the composition
$$
(P^\bullet_B\otimes_AC^\bullet)\otimes_A(P^\bullet_B\otimes_AC^\bullet)
\isom(P^\bullet_B\otimes_AP^\bullet_B)\otimes_AC^\bullet_2
\to P^\bullet_{B_2}\otimes_AC^\bullet\to P^\bullet_B\otimes_AC
$$
where the first isomorphism is again a composition
of associativity isomorphisms and isomorphisms that
swap the factors; the second map is
$\phi_{12}^\bullet\otimes_A\mu^\bullet_C$, where
$\phi_{12}^\bullet:P^\bullet_B\otimes_AP^\bullet_B\to
P^\bullet_{B_2}$ is defined as in \eqref{subsec_pairing-Tors}.
The last map is $P^\bullet_{\mu_B}\otimes_A\one_{C^\bullet}$,
where $P^\bullet_{\mu_B}$ is given by \eqref{subsec_give_Ps}.
Since the associativity and swapping isomorphisms are
obviously natural in all their arguments, we come down
to checking the commutativity, in $\sD(A\Mod)$, of the
diagram of complexes
$$
\xymatrix{
(P^\bullet_B\otimes_AP^\bullet_B)\otimes_AC^\bullet_2
\ar[rrrr]^-{\phi_{12}^\bullet\otimes_A\mu^\bullet_C}
\ar[d]_{(\rho^\bullet_B\otimes_A\rho^\bullet_B)\otimes_A\one_{C^\bullet_2}}
& & & & P^\bullet_{B_2}\otimes_AC^\bullet
\ar[d]^{P^\bullet_{\mu_B}\otimes_A\one_{C^\bullet}} \\
B^\bullet_2\otimes_AC^\bullet_2
\ar[rr]^-{\mu^\bullet_B\otimes_A\mu^\bullet_C} & &
B^\bullet\otimes_AC^\bullet & & P^\bullet_B\otimes_AC^\bullet
\ar[ll]_{\rho^\bullet_B\otimes_A\one_{C^\bullet}}.
}$$
The latter is further reduced to the commutativity of
$$
\xymatrix{
P^\bullet_B\otimes_AP^\bullet_B \ar[rrrr]^-{\phi_{12}^\bullet}
\ar[d]_{\rho^\bullet_B\otimes_A\rho^\bullet_B}
& & & & P^\bullet_{B_2} \ar[d]^{P^\bullet_{\mu_B}} \\
B^\bullet_2 \ar[rr]^-{\mu^\bullet_B} & &
B^\bullet & & P^\bullet_B \ar[ll]_{\rho^\bullet_B}.
}$$
But a simple inspection shows that this diagram indeed
commutes up to homotopy.

Lastly, we claim that if $B^\bullet$ and $C^\bullet$ are
both anti-commutative, then the same holds for
$H_\bullet D^\bullet$. Indeed, consider the diagram
$$
\xymatrix@C-7pt{
(P^\bullet_B\otimes_AP^\bullet_B)\otimes_AC^\bullet_2
\ar[rrrr]^-\sim
\ar[rd]_{(\rho^\bullet_B\otimes_A\rho^\bullet_B)\otimes_A\one_{C^\bullet_2}\ }
\ar@<-4ex>[ddd]^{\phi^\bullet_{12}\otimes_A\mu^\bullet_C}
& & & & (P^\bullet_B\otimes_AP^\bullet_B)\otimes_AC^\bullet_2
\ar[dl]^{\ \ \ (\rho^\bullet_B\otimes_A\rho^\bullet_B)\otimes_A\one_{C^\bullet_2}}
\ar@<4ex>[ddd]_{\phi^\bullet_{12}\otimes_A\mu^\bullet_C} \\
& B_2^\bullet\otimes_AC^\bullet_2 \ar[rr]^-\sim
\ar[rd]_{\mu_B^\bullet\otimes_A\mu^\bullet_C} & &
B_2^\bullet\otimes_AC^\bullet_2 \ar[ld]^{\mu_B^\bullet\otimes_A\mu^\bullet_C} \\
& & B^\bullet\otimes_AC^\bullet \\
P_{B_2}^\bullet\otimes_AC^\bullet
\ar[rr]^-{P^\bullet_{\mu_B}\otimes_A\one_{C^\bullet}} & &
P^\bullet_B\otimes_AC^\bullet
\ar[u]_{\rho^\bullet_B\otimes_A\one_{C^\bullet}} & &
P^\bullet_{B_2}\otimes_AC^\bullet
\ar[ll]_-{P^\bullet_{\mu_B}\otimes_A\one_{C^\bullet}}.
}$$
(The upper isomorphism is obtained from the automorphism
of $P^\bullet_B\otimes_AP^\bullet_B$ that swaps the two factors,
and from the automorphism of $C^\bullet_2$ of the same type;
likewise for the lower isomorphism.)
The assumption on $B^\bullet$ and $C^\bullet$ says that
the inner triangular subdiagram commutes, and the same holds
-- up to homotopy -- for the upper and the left and right
subdiagrams, by a simple inspection. Then also the outer
rectangular subdiagram commutes up to homotopy, and the
contention follows easily.

\begin{remark}\label{rem_from-simpl-to-dga}
(i)\ \
Simplicial $A$-algebras are an important source of
differential graded algebras, thanks to the following
construction. Let $R$ be any simplicial $A$-algebra, $R_\bullet$
the associated chain complex of $A$-modules, and denote
by $\mu_R:R\otimes_AR\to R$ the multiplication map of $R$.
By considering the shuffle map for the bisimplicial $A$-module
$R\boxtimes_AR$ (notation as in \eqref{subsec_shuffle-commute}),
we deduce a natural map of complexes
$$
\mu_{R_\bullet}:
R_\bullet\otimes_AR_\bullet\xrightarrow{\ \Sh_\bullet^{R\boxtimes_AR}\ }
(R\otimes_AR)_\bullet\xrightarrow{\ (\mu_R)_\bullet\ } R_\bullet
$$
and taking into account propositions \ref{prop_commu-shuffle}
and \ref{prop_asso-shuffle}, it is easily seen that
$(R_\bullet,\mu_{R_\bullet})$ is a strictly anti-commutative
differential graded algebra. By remark \ref{rem_dgas}(iii),
we deduce that the graded $A$-module
$H_\bullet R:=\bigoplus_{p\in\N}H_pR$ is naturally an $\N$-graded
associative unital and strictly anti-commutative $A$-algebra.

(ii)\ \
Likewise, if $M$ is any $R$-module, then we obtain on
the associated chain complex $M_\bullet$ a natural
structure of $R_\bullet$-bimodule, so that $H_\bullet M$
is naturally a graded $H_\bullet R$-bimodule.

(iii)\ \
Clearly, a morphism $\phi:R\to S$ of simplicial $A$-algebras
induces a morphism $\phi_\bullet:R_\bullet\to S_\bullet$ of
differential graded $A$-algebras, and a morphism $f:M\to N$
of $R$-modules induces a morphism
$\phi_\bullet:M_\bullet\to N_\bullet$ of $R_\bullet$-bimodules.
\end{remark}

\subsection{Koszul algebras and regular sequences}
\label{sec_koszul-alg}
For any ring $A$, let us consider the category
$$
A\dAlgMod
$$
whose objects are the triples $(B,M,\partial)$, where
$B$ is an $A$-algebra, $M$ is a $B$-module, and
$\partial:M\to B$ is an $A$-linear map. The morphisms
$(B,M,\partial)\to(B',M',\partial')$ in $A\dAlgMod$
are the pairs $(f,g)$ where $f:B\to B'$ is a map of
$A$-algebras, and $g:B'\otimes_BM\to M'$ is a $B'$-linear
map, such that
$$
f\circ\partial(m)=\partial'\circ g(1\otimes m)
\qquad
\text{for every $m\in M$}
$$
with the obvious composition law for such morphisms.
We have a natural functor :
$$
A\dgalt\to A\dAlgMod
\qquad
(B^\bullet,d^\bullet_B)\mapsto(\Ker\,d^0_B,B^{-1},d^{-1}_B)
$$
which admits a left adjoint
$$
\bK_\bullet:A\dAlgMod\to A\dgalt
\qquad
(B,M,\partial)\mapsto\bK_\bullet(B,M,\partial)
$$
that associates to every object $(B,M,\partial)$ its
{\em Koszul algebra} $\bK_\bullet(B,M,\partial)$, whose
underlying strictly anti-commutative $A$-algebra is the
exterior algebra $\Lambda^\bullet_BM$, as in remark
\ref{rem_alternating-alg}(i). The differentials of the
complex $\bK_\bullet(B,M,\partial)$ are given by the maps
$$
d_{n+1}:\Lambda^{n+1}_BM\to\Lambda^n_BM
\qquad
x_0\wedge\cdots\wedge x_n\mapsto\sum_{i=0}^n(-1)^i\partial(x_i)
\cdot x_0\wedge\cdots\wedge x_{i-1}\wedge x_{i+1}\wedge\cdots\wedge x_n
$$
for every $n\in\N$. The detailed verifications shall
be left to the reader.

\begin{remark}\label{rem_koszul-alg}
(i)\ \
In light of remark \ref{rem_alternating-alg}(iii), it is easily
seen that the coproduct of any two objects $(B,M,\partial)$ and
$(B',M',\partial')$ of $A\dAlgMod$ is representable by the object
$$
(B'':=B\otimes_AB',M'':=(M\otimes_AB')\oplus(B\otimes_AM'),\partial'')
$$
where $\partial'':M''\to B''$ is the $A$-linear map
such that
$$
\partial''(m\otimes b',b\otimes m'):=
\partial(m)\otimes b'+b\otimes\partial'(m')
\qquad
\text{for every $m\in M$, $m'\in M'$, $b\in B$, $b'\in B'$}.
$$
Combining with remark \ref{rem_diff-products} and proposition
\ref{prop_was-get-maddd}(iv), we deduce a natural isomorphism :
$$
\bK_\bullet(B'',M'',\partial'')\isom
\bK_\bullet(B,M,\partial)\otimes_A\bK_\bullet(B',M',\partial')
\qquad
\text{in $A\dgalt$}.
$$

(ii)\ \
Let $\bff:=(f_i~|~i=1,\dots,r)$ be a finite system of elements
of $A$, and $I\subset A$ the ideal generated by $\bff$; let
also $\partial_\bff:A^{\oplus r}\to A$ be the $A$-linear form
such that $\partial_\bff(a_1,\dots,a_r):=\sum_{i=1}^rf_ia_i$ for
every $a_1,\dots,a_r\in A$. The {\em Koszul complex of the
sequence $\bff$} is the complex :
$$
\bK_\bullet(\bff):=\bK_\bullet(A,A^{\oplus r},\partial_\bff)
$$
(see \cite[Ch.III, \S1.1]{EGAIII}). Thus, for $r=1$, so
that $\bff=(f)$ for a single element $f\in A$, we have
$$
\bK_\bullet(f)=(0\to A\xrightarrow{\ f\ }A\to 0)
$$
concentrated in homological degrees $0$ and $1$. In the
general case, combining with (i) we get a natural isomorphism :
$$
\bK_\bullet(\bff)\isom
\bK_\bullet(f_1)\otimes_A\cdots\otimes_A\bK_\bullet(f_r)
\qquad
\text{in $A\dgalt$}.
$$
For every complex of $A$-modules $M^\bullet$, we also use
the customary notation :
$$
\bK_\bullet(\bff,M^\bullet):=M^\bullet\otimes_A\bK_\bullet(\bff)
\qquad
\bK^\bullet(\bff,M^\bullet):=
\Tot^\bullet(\Hom_A^\bullet(\bK_\bullet(\bff),M^\bullet))
$$
and denote by $H_\bullet(\bff,M^\bullet)$ (resp.
$H^\bullet(\bff,M^\bullet)$) the homology of
$\bK_\bullet(\bff,M^\bullet)$ (resp. the cohomology
of $\bK^\bullet(\bff,M^\bullet)$). Especially, if $M$
is any $A$-module :
$$
H_0(\bff,M)=M/IM
\qquad
H^0(\bff,M)=\Hom_A(A/I,M)=\Ann_M(I)
$$
(where, as usual, we regard $M$ as a complex placed in degree $0$).
\end{remark}

\begin{lemma}\label{lem_koszul-vanish}
With the notation of remark {\em\ref{rem_koszul-alg}(ii)},
let $g:A\to B$ be a ring homomorphism, and $M^\bullet$
(resp. $N^\bullet$) a complex of $A$-modules (resp.
$B$-modules). The following holds :
\begin{enumerate}
\item
For every $x\in I$, scalar multiplication by $x$ induces
the zero endomorphism on the objects $\bK_\bullet(\bff,M^\bullet)$
and\/ $\bK^\bullet(\bff,M^\bullet)$ of\/ $\Hot(A\Mod)$.
\item
Especially, $H_i(\bff,M^\bullet)$ and\/ $H^i(\bff,M^\bullet)$
are $A/I$-modules, for every $i\in\Z$.
\item
If $I=A$, the complexes $\bK_\bullet(\bff,M^\bullet)$
and\/ $\bK^\bullet(\bff,M^\bullet)$ are homotopically trivial.
\item
Denote by $g(\bff)$ the image in $B$ of the sequence $\bff$.
Then we have natural identifications
$$
\bK_\bullet(\bff,N^\bullet)\isom\bK_\bullet(g(\bff),N^\bullet)
\qquad
\bK^\bullet(\bff,N^\bullet)\isom\bK^\bullet(g(\bff),N^\bullet)
\qquad
\text{in $\sC(B\Mod)$}.
$$
\item
We have a natural isomorphism
$$
\bK^\bullet(\bff,M^\bullet)\isom\bK_\bullet(\bff,M^\bullet)[-r]
\qquad
\text{in $\sC(A\Mod)$}.
$$
\end{enumerate}
\end{lemma}
\begin{proof}(i): It suffices to notice that scalar
multiplication by $f$ induces the zero endomorphism
of $\bK_\bullet(f)$, for every $f\in A$ : indeed, a
homotopy from $f\cdot\one_{\bK_\bullet(f)}$ to the zero
endomorphism is given by the system of maps
$(s_n~|~n\in\N)$ with $s_n:=0$ for $n\neq 0$, and
$s_0:=\one_A$.

(ii) and (iii) are immediate consequences of (i), and
(iv) follows directly from the definitions.

(v): Recall that for every free $A$-module $L$ of
finite rank and every $A$-module $M$ we have a natural
isomorphism of $A$-modules
$$
\Hom_A(L,M)\isom L^\vee\otimes_AM
\qquad
\text{where $L^\vee:=\Hom_A(L,A)$}.
$$
There follows, for every complex $L_\bullet$ of free
$A$-modules of finite rank, a natural isomorphism
$$
\Hom_A^\bullet(L_\bullet,M^\bullet)\isom
L^\vee_\bullet\boxtimes_AM^\bullet
\qquad
\text{in $\sC_2(A\Mod)$}
$$
(notation of example \ref{ex_monoidal}(i)). Moreover, for
any two free $A$-modules of finite rank $L_1$ and $L_2$ we
have as well a natural isomorphism
$$
(L_1\otimes_AL_2)^\vee\isom L_1^\vee\otimes_AL^\vee_2
$$
whence, for any complexes $L_{1,\bullet}$ and $L_{2,\bullet}$
of free $A$-modules of finite rank, a natural isomorphism
$$
(L_{1,\bullet}\boxtimes_AL_{2,\bullet})^\vee\isom
L^\vee_{1,\bullet}\boxtimes_AL^\vee_{2,\bullet}
\qquad
\text{in $\sC_2(A\Mod)$}.
$$
Summing up, it then suffices to observe that there
is a natural isomorphism
$$
\bK_\bullet(f)\isom\bK_\bullet(f)^\vee[1]
\qquad
\text{in $\sC(A\Mod)$}
$$
for every $f\in A$.
\end{proof}

\sset\subsubsection{}\label{subsec_Koszulprime}
Let $\bff$ and $I\subset A$ be as in remark
\ref{rem_koszul-alg}(ii); set as well
$\bff':=(f_1,\dots,f_{r-1})$, and let $I'\subset I$ be the
subideal generated by the sequence $\bff'$. We have a short
exact sequence of complexes :
$0\to A[0]\to\bK_\bullet(f_r)\to A[1]\to 0$,
and after tensoring with $\bK_\bullet(\bff')$, in view of
remark \ref{rem_koszul-alg}(ii) we get a distinguished
triangle :
$$
\bK_\bullet(\bff')\to\bK_\bullet(\bff)\to\bK_\bullet(\bff')[1]
\xrightarrow{\partial}\bK_\bullet(\bff')[1].
$$
By inspecting the definitions one checks easily that the
boundary map $\partial$ is induced by multiplication by
$f_r$. There follow exact sequences :
\set\begin{equation}\label{eq_dist-triang}
0\to H_0(f_r,H_p(\bff',M))\to H_p(\bff,M)\to
H^0(f_r,H_{p-1}(\bff',M))\to 0
\end{equation}
for every $A$-module $M$ and for every $p\in\N$, whence
the following :

\begin{lemma}\label{lem_Kos-cptls-sec}
With the notation of \eqref{subsec_Koszulprime}, the
following conditions are equivalent :
\begin{enumerate}
\alphaenu
\item
$H_i(\bff,M)=0$ for every $i>0$.
\item
The scalar multiplication by $f_r$ is a bijection on
$H_i(\bff',M)$ for every $i>0$, and is an injection on
$M/I'M$.\qed
\end{enumerate}
\end{lemma}

\begin{definition}\label{def_compl-sec}
We say that the sequence $\bff:=(f_1,\dots,f_r)$ of elements
of $A$ is {\em completely secant\/} on the $A$-module $M$, if
we have $H_i(\bff,M)=0$ for every $i>0$.
\end{definition}

The interest of definition \ref{def_compl-sec} is due to its
relation to the notion of {\em regular sequence\/} of elements
of $A$ (see {\em e.g.} \cite[Ch.X, \S9, n.6]{BouAC10}). Namely,
we have the following criterion :

\begin{proposition}\label{prop_Kosz-cptl-sec}
With the notation of \eqref{subsec_Koszulprime}, the following
conditions are equivalent :
\begin{enumerate}
\alphaenu
\item
The sequence $\bff$ is $M$-regular.
\item
For every $j\leq r$, the sequence $(f_1,\dots,f_j)$ is
completely secant on $M$.
\end{enumerate}
\end{proposition}
\begin{proof} Lemma \ref{lem_Kos-cptls-sec} shows that (b) implies
(a). Conversely, suppose that (a) holds; we show that
(b) holds, by induction on $r$. If $r=0$, there is nothing to
prove. Assume that the assertion is already known for all $j<r$.
Since $\bff$ is $M$-regular by assumption, the same holds for the
subsequence $\bff':=(f_1,\dots,f_{r-1})$, and $f_r$ is regular
on $M/(\bff')M$. Hence $H_p(\bff',M)=0$ for every $p>0$, by
inductive assumption. Then lemma \ref{lem_Kos-cptls-sec}
shows that $H_p(\bff,M)=0$ for every $p>0$, as claimed.
\end{proof}

Notice that any permutation of a completely secant sequence
is again completely secant, whereas a permutation of a regular
sequence is not always regular. As an application of the
foregoing, we point out the following :

\begin{corollary}\label{cor_Kosz-cptl-sec}
With the notation of \eqref{subsec_Koszulprime},
the following holds :
\begin{enumerate}
\item
If a sequence $(f,g)$ of elements of $A$ is $M$-regular,
and $M$ is $f$-adically separated, then $(g,f)$ is $M$-regular.
\item
If $(n_1,\dots,n_r)$ is any sequence of strictly positive
integers, then $\bff$ is completely secant on $M$
(resp. $M$-regular) if and only if the same holds
for the sequence $(f_1^{n_1},\dots,f_r^{n_r})$.
\end{enumerate}
\end{corollary}
\begin{proof}(i): According to proposition \ref{prop_Kosz-cptl-sec},
we only need to show that the sequence $(g)$ is completely
secant, {\em i.e.} that $g$ is regular on $M$. Hence, suppose
that $gm=0$ for some $m\in M$; it suffices to show that
$m\in f^nM$ for every $n\in\N$. We argue by induction on $n$.
By assumption $g$ is regular on $M/fM$, hence $m\in fM$, which
shows the claim for $n=1$. Let $n>1$, and suppose we already
know that $m=f^{n-1}m'$ for some $m'\in M$. Hence $0=gf^{n-1}m'$,
so $gm'=0$ and the foregoing case shows that $m'=fm''$ for
some $m''\in\ M$, thus $m=f^nm''$, as required.

(ii): We deal first with the assertion for completely
secant sequences. First, notice that, since every
permutation of a completely secant sequence is still
completely secant, it suffices to show the assertion
for the sequence of integers $(1,1,\dots,1,n_r)$.
However, set $\bff':=(f_1,\dots,f_{r-1})$; by lemma
\ref{lem_Kos-cptls-sec}, $\bff$ is completely secant
if and only if scalar multiplication by $f_r$ is a
bijection on $H_i(\bff',M)$ for every $i>0$, and
an injection on $M/(\bff')M$. But $f_r$ fulfills
the latter conditions if and only if $f_r^{n_r}$ does,
whence the contention.
The assertion for $M$-regular sequences follows
from the foregoing, in view of proposition
\ref{prop_Kosz-cptl-sec}.
\end{proof}

\sset\subsubsection{}\label{subsec_Koszul-is-back}
Let $A$ be a ring, $\bff:=(f_1,\dots,f_r)$ a completely
secant sequence of elements of $A$ (see definition
\ref{def_compl-sec}); denote by $J$ the ideal generated by
$\bff$, and set $A_0:=A/J$. We may regard the complex
$A_0[0]$ as an (especially simple) differential graded
$A$-algebra, with multiplication $\bar\mu{}^\bullet$ deduced
from that of $A_0$, in the obvious way. We may then state :

\begin{proposition}\label{prop_Koszul-is-back}
In the situation of \eqref{subsec_Koszul-is-back}, the
$A_0$-module $J/J^2$ is free of rank $r$, and there exists
a unique isomorphism of strictly anti-commutative graded
$A$-algebras
\set\begin{equation}\label{eq_unique-iso}
\Lambda^\bullet_{A_0}(J/J^2)\isom
H_\bullet(A_0[0]\derotimes_AA_0[0])
\end{equation}
which restricts, in degree $1$, to the natural identification
$$
\Lambda^1_{A_0}(J/J^2)=J/J^2\isom\Tor^A_1(A_0,A_0).
$$
\end{proposition}
\begin{proof} By assumption, the Koszul complex,
with its natural augmentation, yields a resolution
$$
\eps^\bullet:\bK_\bullet(\bff)\to A_0[0]
$$
by free $A$-modules. Hence, there is a unique isomorphism
$\omega^\bullet:P^\bullet_{A_0}\isom\bK_\bullet(\bff)$ in
$\sD(A\Mod)$, whose composition with $\eps^\bullet$ agrees
with $\rho^\bullet_{A_0}$ (notation of \eqref{subsec_give_Ps}).
Then $\eps^\bullet$ is a map of differential graded $A$-algebras,
and we easily deduce a commutative diagram in $\sD(A\Mod)$
$$
\xymatrix{
P^\bullet_{A_0}\otimes_AP^\bullet_{A_0}
\ar[rr]^-{P^\bullet_{\bar\mu}}
\ar[d]_{\omega^\bullet\otimes_A\omega^\bullet} & &
P^\bullet_{A_0} \ar[d]^{\omega^\bullet} \\
\bK_\bullet(\bff)\otimes_A\bK_\bullet(\bff) \ar[rr] & &
\bK_\bullet(\bff)
}$$
whose bottom horizontal arrow is the multiplication map
of $\bK_\bullet(\bff)$, and where $P^\bullet_{\bar\mu}$ is
defined as in \eqref{subsec_give_Ps}. We conclude that
$\omega^\bullet$ induces an isomorphism of anti-commutative
graded $A$-algebras
\set\begin{equation}\label{eq_diarry}
H_\bullet(A_0[0]\derotimes_AA_0[0])\isom\bK_\bullet(\bff,A_0[0])
\isom
H_\bullet((\Lambda^\bullet_A(A^{\oplus r}),d_{\bff,\bullet})\otimes_AA_0).
\end{equation}
By simple inspection, we see that $d_{\bff,i}\otimes_A\one_{A_0}=0$
for every $i\in\Z$, whence an isomorphism of strictly anti-commutative
$A$-graded algebras :
$$
H_\bullet((\Lambda^\bullet_A(A^{\oplus r}),d_{\bff,\bullet})\otimes_AA_0)
\isom
\Lambda^\bullet_{A_0}(A_0^{\oplus r}).
$$
Combining with \eqref{eq_diarry}, we obtain in degree one
a natural isomorphism
$$
\Lambda^1_{A_0}(A_0^{\oplus r})=A_0^{\oplus r}\isom
\Tor^A_1(A_0,A_0)\isom J/J^2
$$
which, finally, delivers the sought isomorphism of
differential graded $A$-algebras. The uniqueness of
\eqref{eq_unique-iso} is clear, since the exterior
algebra is generated by its degree one summand.
\end{proof}

\begin{definition}\label{def_quasi-regular-seq}
Let $A$ be a ring, $\bff:=(f_1,\dots,f_r)$ a finite sequence
of elements of $A$, and $M$ an $A$-module. Set also
$R:=\Z[T_1,\dots,T_r]$, let $I\subset R$ be the ideal
generated by $T_1,\dots,T_r$, and denote by $\gr_\bullet R$
the graded ring associated with the $I$-adic filtration
of $R$; we associate with $\bff$ the $R$-module structure on
$M$ such that $T_ix:=f_ix$ for every $x\in M$ and every
$i=1,\dots,r$, and denote likewise by $\gr_\bullet M$ the graded
$\gr_\bullet R$-module associated with the $I$-filtration on
$M$. We say that $\bff$ is {\em $M$-quasi-regular}, if the
natural map
$$
\gr_\bullet R\otimes_{\gr_0R}\gr_0M\to\gr_\bullet M
$$
is an isomorphism.
\end{definition}

The following result summarizes and completes the list of
interdependencies found thus far between the properties
of finite sequences of elements in a ring that have been
introduced at various stages in the text.

\begin{proposition}
In the situation of definition {\em\ref{def_quasi-regular-seq}},
consider the following conditions :
\begin{enumerate}
\alphaenu
\item
The sequence $\bff$ is $M$-regular.
\item
The sequence $\bff$ is completely secant on $M$.
\item
$\Tor_1^R(R/I,M)=0$.
\item
The sequence $\bff$ is $M$-quasi-regular.
\end{enumerate}
Then {\em (a)$\Rightarrow$(b)$\Rightarrow$(c)$\Rightarrow$(d)},
and if $M$ is $I$-adically complete and separated,
then {\em (d)$\Rightarrow$(a)}.
\end{proposition}
\begin{proof}(a)$\Rightarrow$(b) is already known, by proposition
\ref{prop_Kosz-cptl-sec}.

(b)$\Rightarrow$(c): Let $g:R\to A$ be the unique ring homomorphism
such that $g(T_i):=f_i$ for $i=1,\dots,r$; clearly the $R$-module
structure of $M$ is obtained by restriction of scalars along $g$
from its $A$-module structure. Then, set $\bT:=(T_1,\dots,T_r)$;
by (b) and lemma \ref{lem_koszul-vanish}(iv), we have
$H_i(\bT,M)=0$ for every $i>0$. On the other hand, $\bT$ is
a regular sequence on $R$, hence $H_i(\bT,R)=0$ for every $i>0$,
by the foregoing; {\em i.e.} $\bK_\bullet(\bT)$ is a free resolution
of the $R$-module $R/I$. There follows a natural isomorphism
$H_i(\bT,M)\isom\Tor_i^R(R/I,M)$ for every $i\in\N$, whence
the assertion.

(c)$\Rightarrow$(d): The assumption means that the natural map
$I\otimes_RM\to IM$ is bijective. On the other hand, notice that
$\gr_0:=R/I$ is isomorphic to $\Z$, and $\gr_nR:=I^n/I^{n+1}$ is
a free $\Z$-module for every $n\in\N$; consider then for every
$n\in\N$ the commutative diagram :
$$
\xymatrix{ 0 \ar[r] & I^{n+1}\otimes_RM \ar[r] \ar[d] &
I^n\otimes_RM \ar[r] \ar[d] & \gr_nR\otimes_RM \ar[d] \ar[r] & 0 \\
0 \ar[r] & I^{n+1}M \ar[r] & I^nM \ar[r] & \gr_nM \ar[r] & 0
}$$
all of whose vertical arrows are surjections, and whose
horizontal rows are short exact sequences, as
$\Tor_1^R(\gr_nR,M)=0$. Then, a simple induction shows that
all three vertical arrows are isomorphisms for every $n\in\N$,
whence the assertion.

Suppose now that $M$ is $I$-adically complete and separated,
and that $\bff$ is $M$-quasi-regular; we wish to check that
$\bff$ is $M$-regular, and we shall argue by induction on
the length $r$ of $\bff$. For $r=0$, there is nothing to prove.
Next, suppose that $r\geq 1$, and the assertion is already known
for every $M$-quasi-regular sequence of length $r-1$; we notice

\begin{claim}\label{cl_Chana}
$f_1$ is $M$-regular, and $f_1M\cap I^{n+1}M=f_1I^nM$
for every $n\in\N$.
\end{claim}
\begin{pfclaim} It is easily seen that $(I^{n+1}+T_1I^n)/I^{n+2}$
is a free $\Z$-module for every $n\in\N$, hence we get a short
exact sequence
$$
0\to\gr_nR\otimes_{\gr_0R}\gr_0M\xrightarrow{\phi}
\gr_{n+1}R\otimes_{\gr_0R}\gr_0M\to
((I^{n+1}+T_1I^n)/I^{n+2})\otimes_{\gr_0R}\gr_0M\to 0
$$
for every $n\in\N$ where $\phi$ is induced by scalar
multiplication by $T_1$ on $R$, which maps $I^n$ into $I^{n+1}$
for every such $n$. Then (d) implies that scalar multiplication
by $f_1$ on the $A$-module $M$ induces an injective map
$\gr_nM\to\gr_{n+1}M$ for every $n\in\N$, and since the
$I$-filtration is separated on $M$, we conclude already
that $f_1$ is $M$-regular; moreover :
$$
\{x\in I^nM~|~f_1x\in I^{n+2}M\}=I^{n+1}M
\qquad
\text{for every $n\in\N$}
$$
whence $f_1I^nM\cap I^{n+2}M=f_1I^{n+1}M$ for every $n\in\N$.
By a simple induction on $k$, we deduce that
$f_1I^{n-k}M\cap I^{n+2}M=f_1I^{n+1}M$ for $k=0,\dots,n$, whence
the sought identity.
\end{pfclaim}

Claim \ref{cl_Chana} implies that the $I$-adic topology on $f_1M$
agrees with the topology induced by the $I$-adic topology of $M$;
by virtue of proposition \ref{prop_replaces-Mat-Th.8.1}(i,v), we
deduce that $\bar M:=M/f_1M$ is $I$-adically complete and separated.
Now, let $R':=\Z[T_2,\dots,T_r]\subset R$, set $I':=I\cap R'$, and
denote by $\gr_\bullet R'$ the graded ring associated with the
$I'$-adic filtration of $R'$. The projection $R\to R/(T_1)\isom R'$
induces isomorphisms of rings and respectively graded rings :
$$
\gr_0R\isom\gr_0R'
\qquad
(\gr_\bullet R)/(\bar T_1)\isom\gr_\bullet R'
$$
where $\bar T_1\in\gr_1R$ denotes the class of $T_1$. Let us
endow $\bar M$ with the $R'$-module structure obtained by
restriction of scalars along the inclusion map $R'\to R$,
and let $\gr_\bullet\bar M$ the graded $\gr_\bullet R'$-module
associated with the $I'$-adic filtration on $\bar M$. Then
clearly the $I'$-adic topology of $\bar M$ coincides with
its $I$-adic topology; moreover, notice that
$\gr_\bullet R'\otimes_{\gr_\bullet R}\gr_nM=I^nM/(f_1I^{n-1}M+I^{n+1}M)$,
whereas $\gr_n\bar M=I^nM/((f_1M\cap I^nM)+I^{n+1}M)$ for every
$n>0$. Hence, claim \ref{cl_Chana} also implies that the
projection $\gr_\bullet M\to\gr_\bullet\bar M$ induces isomorphisms of
$\gr_0R$-modules and respectively graded $\gr_\bullet R'$-modules :
$$
\gr_0M\isom\gr_0\bar M
\qquad
\gr_\bullet R'\otimes_{\gr_\bullet R}\gr_\bullet M\isom
\gr_\bullet\bar M
$$
Therefore, the natural map
$$
\gr_\bullet R'\otimes_{\gr_0R'}\gr_0\bar M\to\gr_\bullet\bar M
$$
is also an isomorphism, {\em i.e.} the sequence
$\bff':=(f_2,\dots,f_r)$ is $\bar M$-quasi-regular; by inductive
assumption, $\bff'$ is then $\bar M$-regular, so finally $\bff$
is $M$-regular.
\end{proof}

For future reference, we point out :

\begin{lemma}\label{lem_from-SGA6}
Let $A$ be a ring, $n\geq 1$ an integer, and
$\bff:=(f_0,\dots,f_n)$ a completely secant sequence
of elements of $A$. The following holds :

{\em (i)}\ \
The sequences $\bg:=(f_0,f_1-f_0T_1,\dots,f_n-f_0T_n)$ and
$\bg':=(f_1-f_0T_1,\dots,f_n-f_0T_n)$ are completely secant
in the polynomial $A$-algebra $A[T_1,\dots,T_n]$.

{\em (ii)}\ \
Let $A[f_1/f_0,\dots,f_n/f_0]\subset A[1/f_0]$ be the
$A$-subalgebra generated by $f_1/f_0,\dots,f_n/f_0$. Then
the kernel of the surjective map of $A$-algebras
$$
A[T_1,\dots,T_n]\to A[f_1/f_0,\dots,f_n/f_0]
\qquad
T_i\mapsto f_i/f_0
\qquad\text{for $i=1,\dots,n$}
$$
is the ideal $I\subset A[T_1,\dots,T_n]$ generated by
the sequence $\bg'$.
\end{lemma}
\begin{proof}(i): Since the inclusion map
$A\to B:=A[T_1,\dots,T_n]$ is flat, it is clear that the
sequence $\bff$ is completely secant in $B$. Define the
$B$-linear forms
$\partial_\bff,\partial_\bg:B^{\oplus n+1}\to B$ as in
remark \ref{rem_koszul-alg}(ii), so that
$\bK_\bullet(\bff)=\bK_\bullet(B,B^{\oplus n+1},\partial_\bff)$,
and likewise for $\bK_\bullet(\bg)$. We have an isomorphism
$$
(\one_B,\phi):(B,B^{\oplus n+1},\partial_\bg)\isom
(B,B^{\oplus n+1},\partial_\bff)
\qquad
\text{in $A\dAlgMod$}
$$
with $\phi:B^{\oplus n+1}\isom B^{\oplus n+1}$ the $B$-linear
automorphism such that $\phi(e_0):=e_0$ and
$\phi(e_i):=e_i-T_ie_0$ for $i=1,\dots,n$. The assertion
for $\bg$ is an immediate consequence.

Combining with lemma \ref{lem_Kos-cptls-sec}, we see that
the scalar multiplication by $f_0$ is injective on
$H_i(\bg',B)$, for every $i\in\N$. Hence, the natural map
$$
H_i(\bg',B)\to H_i(\bg',B[f_0^{-1}])
$$
is injective for every $i\in\N$, and in order to show that
$\bg'$ is completely secant in $B$, we may then assume that
$f_0\in A^\times$. To this aim, we consider also the sequence
$\bT:=(T_1,\dots,T_n)$; we have an isomorphism
$$
(\psi,\one_{B^{\oplus n}}):(B,B^{\oplus n},\partial_\bT)
\isom(B,B^{\oplus n},\partial_{\bg'})
\qquad
\text{in $A\dAlgMod$}
$$
where $\psi:B\isom B$ is the ring isomorphism such that
$\psi(T_i):=f_i-f_0T_i$ for $i=1,\dots,n$. Since the sequence
$\bT$ is obviously $B$-regular, the assertion then follows
from proposition \ref{prop_Kosz-cptl-sec}.

(ii)\ \
We have a commutative diagram of rings :
$$
\xymatrix{ B/IB \ar[r] \ar[d] & A[f_1/f_0,\dots,f_n/f_0] \ar[d] \\
B/IB[1/f_0] \ar[r] & A[1/f_0]
}$$
whose vertical arrows are the localizations, and it is
easily seen that the bottom horizontal arrow is an isomorphism.
On the other hand, from (i) and lemma \ref{lem_Kos-cptls-sec}
we see that the scalar multiplication by $f_0$ is injective
on $B/I$; hence the top horizontal arrow is injective, whence
the contention.
\end{proof}

\begin{definition}\label{def_essentially-zero}
Let $\cA$ be any additive category, and $A_\bullet$ any
object of $\bFun(\N^o,\cA)$ (where the partially ordered
set $(\N,\leq)$ is regarded as category, as in example
\ref{ex_universe}(iii)). In other words, $A_\bullet$ is a
datum $(A_\bullet,\phi_{\bullet\bullet})$ consisting of a family
$(A_p~|~p\in\N)$ of objects of $\cA$ and a system of
morphisms $\phi_{q,p}:A_q\to A_p$ of $\cA$ for every
integers $q\geq p\geq 0$, such that
$$
\phi_{q,p}\circ\phi_{r,q}=\phi_{r,p}
\qquad
\text{for every $r\geq q\geq p\geq 0$}.
$$
The objects of $\bFun(\N^o,\cA)$ are called also
{\em inverse systems} of $\cA$ indexed by $\N$. Then :

(i)\ \
We say that $A_\bullet$ is {\em essentially zero}, if for
every $p\in\N$ there exists $q\geq p$ such that $\phi_{q,p}$
is the zero morphism.

(ii)\ \
We say that $A_\bullet$ is {\em uniformly essentially zero},
if there exists $c\in\N$ such that $\phi_{p+c,p}$ is the
zero morphism for every $p\in\N$. In this case, the smallest
of such $c$ is called the {\em step} of $A_\bullet$.
\end{definition}

Recall that if $\bFun(\N^o,\cA)$ is an additive category,
and it is even abelian, if the same holds for $\cA$; moreover,
the kernels and cokernels in $\bFun(\N^o,\cA)$ are computed
argumentwise : see remark \ref{rem_Add.Fun}(ii). We notice :

\begin{lemma}\label{lem_inverse-Serre-subcat}
Let $\cA$ be any abelian category, and consider a short
exact sequence
$$
0\to A'_\bullet\to A_\bullet\to A''_\bullet\to 0
\qquad
\text{in $\bFun(\N^o,\cA)$}.
$$
Then, $A_\bullet$ is essentially zero (resp. uniformly essentially
zero) if and only if the same holds for both $A'_\bullet$ and
$A''_\bullet$.
\end{lemma}
\begin{proof} Indeed, it is clear that if $A_n\to A_m$ is
the zero morphism for some $n\geq m$, then the same holds
for the morphisms $A'_n\to A'_m$ and $A''_n\to A''_m$.
Conversely, suppose that $A'_n\to A'_m$ and $A''_p\to A''_n$
are the zero morphisms, for some $p\geq n\geq m$; then it
follows easily that the same holds for the morphism
$A_p\to A_m$.
\end{proof}

\sset\subsubsection{}\label{subsec_def-bphis}
Let $\bff:=(f_1,\dots,f_r)$ and $\bg:=(g_i~|~i=1,\dots,r)$ be
two finite sequences of elements of a ring $A$; we set
$\bff\bg:=(f_ig_i~|~i=1,\dots,r)$ and define as follows
a map of Koszul complexes
$$
\bphi_\bg:\bK_\bullet(\bff\bg)\to\bK_\bullet(\bff)
$$
(see remark \ref{rem_koszul-alg}(ii)). First, suppose that
$r=1$; then $\bff=(f)$, $\bg=(g)$ and the sought map $\bphi_g$
is the commutative diagram :
$$
\xymatrix{
0 \ar[r] & A \ar[r]^-{fg} \ar[d]_g & A \ar[r] \ddouble & 0 \\
0 \ar[r] & A \ar[r]^-f & A \ar[r] & 0.
}$$
For the general case we let :
$$
\bphi_\bg:=
\bphi_{g_1}\otimes_A\cdots\otimes_A\boldsymbol{\phi}_{g_r}.
$$
Especially, for every $m,n\geq 0$ we have maps
$\bphi_{\bff^n}:\bK_\bullet(\bff^{n+m})\to\bK_\bullet(\bff^m)$,
whence maps
$$
\bphi_{\bff^n}^\bullet:
\bK^\bullet(\bff^m,M)\to\bK^\bullet(\bff^{m+n},M)
$$
and clearly $\bphi^\bullet_{\bff^{p+q}}=
\bphi^\bullet_{\bff^p}\circ\bphi^\bullet_{\bff^q}$ for every
$m,p,q\geq 0$.

\sset\subsubsection{}\label{subsec_badabum}
In the situation of \eqref{subsec_def-bphis}, set
$A_r:=\Z[T_1,\dots,T_r]$, and denote by $\beta_\bff:A_r\to A$
the unique ring homomorphism such that $\beta_\bff(T_i):=f_i$
for every $i=1,\dots,r$. Let also $I_r\subset A_r$ be the
ideal generated by $(T_1,\dots,T_r)$. We consider the
following conditions :
\begin{itemize}
\item[$\mathrm{(a)}_\bff$]
The inverse system $(\Tor^{A_r}_i(A_r/I_r^n,A)~|~n\in\N)$
is essentially zero for every $i>0$.
\item[$\mathrm{(b)}_\bff$]
The inverse system $(\Tor^A_i(A/I^n,A/I)~|~n\in\N)$ is
essentially zero for every $i>0$.
\item[$\mathrm{(c)}_\bff$]
The inverse system $(H_i(\bff,I^n)~|~n\in\N)$
is essentially zero for every $i>0$.
\item[$\mathrm{(d)}_\bff$]
The inverse system $(H_i\bK_\bullet(\bff^n)~|~n\in\N)$
is essentially zero for every $i>0$.
\item[$\mathrm{(e)}_\bff$]
The inverse system $(\Tor_{i+1}^A(A/I^n,M)~|~n\in\N)$ is
essentially zero, for every $i,k\in\N$ and every
$A/I^k$-module $M$.
\item[$\mathrm{(f)}_\bff$]
For every $i\in\Z$ and every bounded above complex $C_\bullet$
of flat $A$-modules such that for every $j\in\Z$ the
annihilator ideal of $H_jC_\bullet=0$ contains a power
of $I$, the inverse system $(H_i(I^nC_\bullet)~|~n\in\N)$
is essentially zero.
\end{itemize}
Moreover, we shall also consider the conditions
$\mathrm{(a)^{un}_\bff,\dots,(f)^{un}_\bff}$ obtained
from $\mathrm{(a)_\bff,\dots,(f)_\bff}$ after replacing
``essentially zero'' by ``uniformly essentially zero''.

\begin{remark}
Condition $\mathrm{(d)_\bff}$ seems to have been first
considered in \cite[Lemma 2.4]{LCoh}, and in particular,
our remark \ref{rem_Hartsho} was already observed in
\cite[Lemma 2.5 and th.2.8]{LCoh}. The same condition
reappears in \cite[Exp.II, Lemme 9]{SGA2}, as well as in
in \cite[Lemma 3.1.1]{AJL}, and it is baptised
``weak proregularity'' in an {\em errata}\footnote{see
\url{http://www.math.purdue.edu/~lipman/papers/homologyfix.pdf}}
for that paper. A non-commutative extension is proposed in \cite{Yeku}.
For motivation, we point out the following lemma, which says
that completely secant sequences $\bff$ trivially satisfy
condition $\mathrm{(a)^{un}_\bff}$ :
\end{remark}

\begin{lemma}\label{lem_secant-badabum}
With the notation of \eqref{subsec_badabum}, if the sequence
$\bff$ is completely secant, we have
$$
\Tor_i^{A_r}(A_r/I_r^n,A)=0
\qquad
\text{for every $n\in\N$ and every $i>0$}.
$$
\end{lemma}
\begin{proof} Notice first that $A_r/I_r=\Z$, and
$I_r^n/I_r^{n+1}$ is a free $\Z$-module, for every
$n\in\N$; using the long $\Tor$-exact sequences
arising from the short exact sequences
$$
0\to I^n_r/I^{n+1}_r\to A_r/I^{n+1}_r\to A_r/I^n_r\to 0
$$
a simple induction reduces to the case where $n=1$
(details left to the reader). In this case, we consider
the sequence $\bt:=(T_1,\dots,T_r)$ of elements of $A_r$,
which is clearly completely secant, so that the Koszul
complex $\bK_\bullet(\bt)$ is a resolution of the
$A_r$-module $A_r/I_r$ by free $A_r$-modules, whence
natural isomorphisms
$$
\Tor_i^{A_r}(A_r/I_r,A)\isom H_i(\bt,A)\isom H_i(\bff,A)
\qquad
\text{for every $i\in\N$}
$$
(lemma \ref{lem_koszul-vanish}(iv)), and the contention ensues.
\end{proof}

\begin{lemma}\label{lem_ess-zero-Tors}
With the notation of \eqref{subsec_badabum}, for every
$i,p\in\N$ with $i>0$, and every $A_r/I^p_r$-module $N$,
the natural morphism
$\Tor_i^{A_r}(A_r/I_r^{n+p},N)\to\Tor_i^{A_r}(A_r/I_r^n,N)$
is the zero map.
\end{lemma}
\begin{proof} Let $\Fil^\bullet N$ be the $I_r$-adic
filtration of $N$, and denote by $\gr^\bullet N$ the
associated graded $A_r/I_r$-module; using the long exact
$\Tor$-sequences arising from the short exact
sequences $0\to\Fil^{q+1}N\to\Fil^qN\to\gr^qN\to 0$
(for every $q\in\N$), and an easy induction argument,
we reduce to the case where $p=1$ (details left to
the reader). In this case, $N$ is actually an
$A_r/I_r$-module, {\em i.e.} a $\Z$-module; we may
moreover reduce to the case where $N$ is a finitely
generated $\Z$-module, and then we may even assume
that $N=\Z/a\Z$ for some $a\in\N$ (details left to
the reader). In this situation, set $B:=A_r/aA_r$,
$J:=I_rB_r$, consider the standard $2$-spectral
sequence (\cite[Th.5.6.6]{We})
$$
E^2_{pq}:=\Tor^B_p(\Tor^{A_r}_q(A_r/I^n_r,B),\Z/a\Z)
\Rightarrow\Tor^{A_r}_{p+q}(A_r/I^n_r,\Z/a\Z)
$$
and notice that, on the one hand, the scalar multiplication
by $a$ is injective on $A_r/I^n_r$, and on the other hand,
the complex $0\to A_r\to A_r\to 0$ with differential given
by $a\cdot\one_{A_r}$ is a free resolution of $B$; we deduce
that $E^2_{pq}=0$ whenever $q>0$, whence natural isomorphisms :
$$
\Tor^B_i(B/J^n,\Z/a\Z)\isom\Tor^{A_r}_i(A_r/I^n_r,\Z/a\Z)
\qquad
\text{for every $i,n\in\N$}
$$
and we are reduced to checking that the natural map
$$
\Tor^B_i(B/J^{n+1},\Z/a\Z)\to
\Tor^B_i(B/J^n,\Z/a\Z)
$$
vanishes for every $i>0$. To this aim, notice that
the Koszul complex $\bK_\bullet:=\bK_\bullet(T_\bullet,B)$
associated with the sequence $T_\bullet:=(T_1,\dots,T_k)$
of elements of $B$, is a resolution of $Z/a\Z$ consisting
of free $B$-modules (proposition \ref{prop_Kosz-cptl-sec}),
so we have a natural isomorphism
$$
B/J^n\otimes_B\bK_\bullet\isom
B/J^n\derotimes_B\Z/a\Z
\qquad
\text{in $\sD(B\Mod)$}.
$$
Now, denote by $\gr_\bullet B$ the standard graded
$\Z/a\Z$-algebra structure on $B$ such that
$\gr_1B$ is generated by $T_1,\dots,T_k$.
For every $i\in\N$ we define a grading on $\bK_i$ as
follows. Let $(e_1,\dots,e_k)$ be the canonical basis
of $B^{\oplus k}$, and $\cF_i$ the set of all strictly
increasing maps $\{1,\dots,i\}\to\{1,\dots,r\}$; for
every $\phi\in\cF_i$ and every $r\in\N$ set
$$
e_\phi:=e_{\phi(1)}\wedge\cdots\wedge e_{\phi(i)}
\qquad\text{and}\qquad
\gr_r\bK_i:=\sum_{\phi\in\cF_i}e_\phi\cdot\gr_{r-i}B.
$$
A simple inspection then shows that the differential
of the Koszul complex is, in each degree, a map of
graded $\Z/a\Z$-modules, therefore we get a natural
decomposition
$$
\bK_\bullet\isom\bigoplus_{r\in\N}\gr_r\bK_\bullet
\qquad
\text{in $\sC(\Z/a\Z\Mod)$}
$$
and since $J^n=\bigoplus_{r\geq n}\gr_rB$, the
complex $\bK_\bullet\otimes_BB/J^n$ admits
a corresponding decomposition, and we have
$$
(\gr_r\bK_\bullet\otimes_BB/J^n)_i=
\left\{\begin{array}{ll}
         \gr_r\bK_i & \text{if $r<i+n$} \\
         0          & \text{otherwise}.
       \end{array}\right.
$$
Since $H_i\bK_\bullet=0$ for every $i>0$, it follows
easily that
$$
\gr_rH_i(\bK_\bullet\otimes_BB/J^n)=0
\qquad
\text{for every $i,r,n\in\N$ such that $r\neq i-1+n$ and $i>0$}.
$$
Thus, fix $i,n\in\N$ with $i>0$, and consider any
$x\in\Tor_i^B(B/J^{n+1},\Z/a\Z)$; by the foregoing,
we have $x\in H_i(\gr_{i+n}\bK_\bullet\otimes_BB/J^{n+1})$,
and $H_i(\gr_{i+n}\bK_\bullet\otimes_BB/J^n)=0$,
so the image of $x$ vanishes in $\Tor_i^B(B/J^n,\Z/a\Z)$,
as required.
\end{proof}

\begin{proposition}\label{prop_long-equivalences}
With the notation of \eqref{subsec_badabum}, fix $i\in\N$,
and let also $\bg:=(g_i~|~i=1,\dots,s)$ be another sequence
of elements of $A$ that generates an ideal $J$. We have :
\begin{enumerate}
\item
$\mathrm{(a)_\bff\Leftrightarrow(b)_\bff\Leftrightarrow(c)_\bff
\Leftrightarrow(d)_\bff\Leftrightarrow(e)_\bff
\Leftrightarrow(f)_\bff}$.
\item
$\mathrm{(a)^{un}_\bff\Rightarrow(b)^{un}_\bff
\Leftrightarrow(c)^{un}_\bff\Leftrightarrow(e)^{un}_\bff
\Leftrightarrow(f)^{un}_\bff}$.
\item
If the radical of $I$ equals the radical of $J$, then
$\mathrm{(a)_\bff\Leftrightarrow(a)_\bg}$.
\item
If $I=J$, then $\mathrm{(a)^{un}_\bff\Leftrightarrow(a)^{un}_\bg}$.
\end{enumerate}
\end{proposition}
\begin{proof} Let us check first that
$\mathrm{(a)_\bff\Rightarrow(b)_\bff}$ and
$\mathrm{(a)^{un}_\bff\Rightarrow(b)^{un}_\bff}$. To this aim,
we consider the change of ring spectral sequence
$$
E(n)^2_{ij}:=\Tor^A_i(\Tor^{A_r}_j(A_r/I_r^n,A),A/I)\Rightarrow
\Tor_{i+j}^{A_r}(A_r/I_r^n,A/I).
$$
Notice that the projection $A_r/I_r^m\to A_r/I_r^n$ induces
a morphism of spectral sequences
\set\begin{equation}\label{eq_several-planes}
E(m)^\bullet_{\bullet\bullet}\to E(n)^\bullet_{\bullet\bullet}
\qquad
\text{for every $m\geq n$}.
\end{equation}
Then the assertion is the case $p=0$ of the following :

\begin{claim} If condition $\mathrm{(a)}_\bff$ (resp.
$\mathrm{(a)^{un}_\bff}$) holds, then the inverse system
$(E(n)^{p+2}_{i+1,0}~|~n\in\N)$ is essentially zero (resp.
uniformly essentially zero) for every $i,p\in\N$.
\end{claim}
\begin{pfclaim} We fix $i\in\N$, and we argue by
descending induction on $p$. Notice first that
$$
E(n)^{i+1}_{i,0}=E(n)^\infty_{i,0}
\qquad
\text{for every $i,n\in\N$}.
$$
Then, lemma \ref{lem_ess-zero-Tors} trivially implies that the
assertion holds for every $i,p\in\N$ with $p\geq i$. Now, let
$0\leq q\leq i-1$, and suppose that the assertion is already
known for every $p>q$ (with our fixed $i$). We get an exact
sequence of inverse systems
$$
0\to(E(n)^{q+3}_{i+1,0}~|~n\in\N)\to(E(n)^{q+2}_{i+1,0}~|~n\in\N)
\to(E(n)^{q+2}_{i-q-1,q+1}~|~n\in\N)
$$
and by inductive assumption the system
$(E(n)^{q+3}_{i+1,0}~|~n\in\N)$ is essentially zero (resp.
uniformly essentially zero). Moreover, $\mathrm{(a)}_\bff$
(resp. $\mathrm{(a)^{un}_\bff}$) implies that the same holds
for the system $(E(n)^{q+2}_{i-q-1,q+1}~|~n\in\N)$. Then the
contention follows from lemma \ref{lem_inverse-Serre-subcat}.
\end{pfclaim}

$\bullet$\ \
In order to show that $\mathrm{(b)_\bff\Rightarrow(a)_\bff}$
we remark :

\begin{claim}\label{cl_b-implies-a}
Let $j\in\N$ be any integer. Then :
\begin{enumerate}
\item
If $(E(n)^2_{0,j}~|~n\in\N)$ is an essentially zero inverse
system, the same holds for the inverse system
$(\Tor^{A_r}_j(A_r/I^n_r,A)~|~n\in\N)$.
\item
If $(\Tor^{A_r}_j(A_r/I^n_r,A)~|~n\in\N)$ is an essentially
zero inverse system, the same holds for the inverse system
$(E^{p+2}_{i,j}~|~n\in\N)$, for every $i,p\in\N$.
\item
The inverse system $(E(n)^{j+2}_{0,j}~|~n\in\N)$ is essentially
zero.
\item
$\mathrm{(b)}_\bff$ implies that $(E(n)^{j+2}_{i,0}~|~n\in\N)$
is an essentially zero system for every $i\in\N$.
\end{enumerate}
\end{claim}
\begin{pfclaim} (i): Denote by $\phi_{n,c}:
\Tor_j^{A_r}(A_r/I_r^{n+c},A)\to\Tor^{A_r}_j(A_r/I^n_r,A)$ the
natural map. The hypothesis means that for every $n\in\N$
there exists $c\in\N$ such that $\phi_{n,c}\otimes_AA/I=0$,
{\em i.e.}
$\Img\,\phi_{n,c}\subset I\cdot\Tor^{A_r}_j(A_r/I^n_r,A)$.
By a simple induction we deduce that for every $n,k\in\N$
there exists $d\in\N$ such that
$\Img\,\phi_{n,c}\subset I^k\cdot\Tor^{A_r}_j(A_r/I^n_r,A)$,
and letting $k=n$, the assertion follows.

(ii): Obviously, the assumption implies that
$(E^2_{i,j}~|~n\in\N)$ is an essentially zero inverse system;
the assertion is an immediate consequence.

(iv) is clear, and notice that $E(n)^{j+2}_{0,j}=E(n)^\infty_{0,j}$
for every $n\in\N$. In light of lemma \ref{lem_ess-zero-Tors},
assertion (iii) follows easily as well.
\end{pfclaim}

In light of claim \ref{cl_b-implies-a}(i) it suffices to
show that $(E(n)^2_{0,j})$ is essentially zero for every
$j>0$. We argue by induction on $j$; for $j=1$ we consider
the exact sequence of inverse systems :
$$
(E(n)^2_{2,0}~|~n\in\N)\to(E(n)^2_{0,1}~|~n\in\N)\to
(E(n)^3_{0,1})\to 0.
$$
By claim \ref{cl_b-implies-a}(iii,iv), the first and third
terms are essentially zero, and therefore the same holds
for the middle term, by lemma \ref{lem_ess-zero}.

Next, let $j>1$, and suppose that the inverse systems
$(E(n)^2_{0,k}~|~n\in\N)$ are essentially zero for
$0<k<j$. By claim \ref{cl_b-implies-a}(i,ii) we deduce
that $(E(n)^{p+2}_{i,k}~|~n\in\N)$ is essentially zero
for every $p,i\in\N$ and whenever $0<k<j$. We show, by
descending induction on $p$, that $(E(n)^{p+2}_{0,j}~|~n\in\N)$
is essentially zero for every $p\in\N$. First, the assertion
holds for $p\geq j$, by virtue of claim \ref{cl_b-implies-a}(iii).
Suppose then that $0\leq k<j$, and that the assertion is already
known for every $p>k$; we consider the exact sequence of
inverse systems
$$
(E(n)^{k+2}_{k,j-k-1}~|~n\in\N)\to(E(n)^{k+2}_{0,j}~|~n\in\N)\to
(E(n)^{k+3}_{0,j})\to 0.
$$
Our inductive assumtions imply that the first and third
terms are essentially zero; then the same holds for the
middle term (lemma \ref{lem_ess-zero}), and the proof is
concluded.

$\bullet$\ \
Next, we show that $\mathrm{(b)_\bff\Rightarrow(e)_\bff}$
and $\mathrm{(b)^{un}_\bff\Rightarrow(e)^{un}_\bff}$.
For any $A/I^k$-module $M$, we need to check that
the inverse system $(\Tor_{i+1}^A(A/I^n,M)~|~n\in\N)$
is essentially zero (resp. uniformly essentially zero);
to this aim, we consider the $I$-adic filtration
$\Fil^\bullet M$ on $M$, and the associated graded
$A/I$-module $\gr_\bullet M$. Applying the long exact
$\Tor$-sequences arising from the short exact sequences
$0\to\Fil^{q+1}M\to\Fil^qM\to\gr^qM\to 0$, together with
lemma \ref{lem_inverse-Serre-subcat}, we reduce to checking
that  $(\Tor_{i+1}^A(A/I^n,\gr^\bullet M)~|~n\in\N)$ is an
essentially zero (resp. uniformly essentially zero)
inverse system. Thus, we may suppose that $IM=0$, in
which case we show more precisely :

\begin{claim}\label{cl_precision}
Let $p>0$ be any integer, and suppose that
$(\Tor^A_i(A/I^n,A/I)~|~n\in\N)$ is an essentially
zero (resp. uniformly essentially zero) inverse system for
every $i=1,\dots,p$. Then, for every $A/I$-module $M$
the inverse system $(\Tor^A_p(A/I^n,M)~|~n\in\N)$ is
essentially zero (resp. uniformly essentially zero).
\end{claim}
\begin{pfclaim} We consider the spectral sequence
$$
E(n)^2_{i,j}:=\Tor_i^{A/I}(\Tor^A_j(A/I^n,A/I),M)
\Rightarrow\Tor^A_{i+j}(A/I^n,M)
\qquad
\text{for every $n\in\N$}
$$
and the corresponding morphisms of spectral sequences as
in \eqref{eq_several-planes}, induced by the projections
$A/I^m\to A/I^n$. Notice that $E(n)^2_{i,0}=0$, and
therefore $E(n)^\infty_{i,0}=0$ for every $i>0$. On the other
hand, our assumption implies that the inverse system
$(E(n)^2_{i,j}~|~n\in\N)$ is essentially zero (resp. uniformly
essentially zero) for every $i,j\in\N$ such that $j>0$ and
$i+j=q$. Summing up, we deduce that the inverse system
$(E(n)^\infty_{i,j}~|~n\in\N)$ is essentially zero (resp.
uniformly essentially zero) for every $i,j\in\N$ such that
$i+j=p$. Taking into account lemma \ref{lem_inverse-Serre-subcat},
the contention now follows by a simple induction.
\end{pfclaim}

$\bullet$\ \
We check next that
$\mathrm{(e)}_\bff\Rightarrow\mathrm{(f)}_\bff$ and
$\mathrm{(e)}^{un}_\bff\Rightarrow\mathrm{(f)}^{un}_\bff$.
To this aim, notice that
$H_i(I^nC_\bullet)=H_i(I^n[0]\derotimes_AC_\bullet)$ for every
$i,n\in\Z$, since $C_\bullet$ is a complex of flat $A$-modules.
On the other hand, we have a standard spectral sequence
$$
E(n)^2_{p,q}:=\Tor_p^A(I^n,H_qC_\bullet)\Rightarrow
H_{p+q}(I^n[0]\derotimes_AC_\bullet).
$$
Since
$\Tor_p^A(I^n,H_qC_\bullet)\isom\Tor_{p+1}^A(A/I^n,H_qC_\bullet)$
for every $p,n\in\N$, condition $\mathrm{(e)}_\bff$ (resp.
$\mathrm{(e)^{un}_\bff}$) implies that $(E(n)^2_{p,q}~|~n\in\N)$
is an essentially zero (resp. uniformly essentially zero)
inverse system for every $p\in\N$ and every $q\in\Z$, hence
the same holds for the system $(E(n)^\infty_{p,q}~|~n\in\N)$. Since
$C_\bullet$ is bounded above, lemma \ref{lem_inverse-Serre-subcat}
easily implies that $H_i(I^n[0]\derotimes_AC_\bullet)$ is an
essentially zero (resp. uniformly essentially zero) system
for every $i\in\Z$ (details left to the reader), whence the
contention.

Clearly $\mathrm{(f)}_\bff\Rightarrow\mathrm{(c)}_\bff$ and
$\mathrm{(f)^{un}_\bff\Rightarrow (c)^{un}_\bff}$,
since the Koszul complex is a bounded complex of flat
$A$-modules whose cohomology is an $A/I$-module in each
degree (lemma \ref{lem_koszul-vanish}(ii)).

$\bullet$\ \
To show that $\mathrm{(c)}_\bff\Rightarrow\mathrm{(b)}_\bff$
and $\mathrm{(c)^{un}_\bff\Rightarrow (b)^{un}_\bff}$
we consider the spectral sequence
$$
E(n)^2_{i,j}:=\Tor^A_i(I^n,H_j\bK_\bullet(\bff))
\Rightarrow H_{i+j}(\bff,I^n)
\qquad
\text{for every $n\in\N$}
$$
and the corresponding system of morphisms of spectral sequences
as in \eqref{eq_several-planes}. Recalling that
$H_0\bK_\bullet(\bff_\bullet)=A/I$, as well as the natural isomorphism
\set\begin{equation}\label{eq_yabadaba}
\Tor^A_i(I^n,A/I)\isom\Tor^A_{i+1}(A/I^n,A/I)
\qquad
\text{for every $i,n\in\N$}
\end{equation}
the assertion will follow from the case $j=p=0$
of the following more general :

\begin{claim}
If $\mathrm{(c)}_\bff$ (resp. $\mathrm{(c)^{un}_\bff}$)
holds, the inverse system $(E(n)^{p+2}_{i,j}~|~n\in\N)$ is
essentially zero (resp. uniformly essentially zero) for
every $i,j,p\in\N$.
\end{claim}
\begin{pfclaim} We argue by induction on $i$. For $i=0$,
notice that $E(n)^2_{0,0}=E(n)^\infty_{0,0}$ for every $n\in\N$;
condition $\mathrm{(c)}_\bff$ (resp. $\mathrm{(c)^{un}_\bff}$)
then yields the claim for the system $(E(n)^2_{0,0}~|~n\in\N)$.
In light of \eqref{eq_yabadaba} and claim \ref{cl_precision},
we deduce that $(E(n)^2_{0,j}~|~n\in\N)$ is essentially
zero (resp. uniformly essentially zero) for every $j\in\N$,
and then the same follows for $(E(n)^{p+2}_{0,j}~|~n\in\N)$,
for every $j,p\in\N$.

Next, suppose that $k>0$, and the claim is already known
for every $i<k$ and every $j,p\in\N$. We show, by descending
induction on $q$, that $(E(n)^q_{k,0}~|~n\in\N)$ is an
essentially zero (resp. uniformy essentially zero)
system for every $q\geq 2$. Indeed, notice first that
$E(n)^{k+1}_{k,0}=E(n)^\infty_{k,0}$; in light of
$\mathrm{(c)}_\bff$ (resp. $\mathrm{(c)^{un}_\bff}$), the
assertion follows already for every $q\geq k+1$. Next,
let $p\leq k$ and suppose that the assertion is already
known for $q=p+1$. Notice the exact sequence of inverse
systems :
$$
0\to(E(n)^{p+1}_{k,0}~|~n\in\N)\to
(E(n)^p_{k,0}~|~n\in\N)\to(E(n)^p_{k-p,p-1}~|~n\in\N).
$$
By inductive assumption, both the first and third terms
are essentially zero (resp. uniformly essentially zero),
so the same holds for the middle term. We know now that
$(E(n)^2_{k,0}~|~n\in\N)$ is essentially zero (resp.
uniformly essentially zero); from \eqref{eq_yabadaba} and
claim \ref{cl_precision} we deduce that the same holds for
$(E(n)^2_{k,j}~|~n\in\N)$, for every $j\in\N$, and finally
the same follows for $(E(n)^{p+2}_{k,j}~|~n\in\N)$, for
every $p,j\in\N$, whence the claim.
\end{pfclaim}

$\bullet$\ \
Lastly, let us check that $\mathrm{(a)_\bff\Leftrightarrow(d)_\bff}$.
To this aim, set $T_\bullet:=(T_1,\dots,T_r)$, and notice that
$H_i(T^n_\bullet,A_r)=0$ for every $i,n>0$ (proposition
\ref{prop_Kosz-cptl-sec}), whence a natural isomorphism
$$
\Tor^{A_r}_i(A_r/T^n_\bullet A_r,A)\isom
H_i(\bK_\bullet(T^n_\bullet)\otimes_{A_r}A)\isom H_i(\bff^n,A)
\qquad
\text{for every $i,n\in\N$}.
$$
On the other hand, it is clear that for every $n\in\N$
there exists $m\in\N$ such that
$$
I^m_r\subset T_\bullet^nA_r\subset I^n_r
$$
whence the contention.

(iv): Set $B_r:=A_r\otimes_\Z A$, and notice that $\beta_\bff$
factors as a composition
$$
A_r\xrightarrow{\ \beta'_\bff\ }B_r\xrightarrow{\ \gamma_\bff\ }A_r
$$
where $\beta'_\bff:A_r\to B_r$ is the ring homomorphism such that
$\beta'_\bff(T_i):=1\otimes f_i$ for $i=1,\dots,r$, and
$\gamma_\bff$ is the map of $\Z$-algebras induced by the pair
$(\beta_\bff,\one_A)$. Let us endow $B_r$ with the $A_r$-module
structure given by $\beta'_\bff$ and $A$ with the $B_r$-module
structure given by $\gamma_\bff$; there follows, for every
$n\in\N$ a change of ring spectral sequence
$$
E(n)^2_{i,j}:=\Tor_i^{B_r}(\Tor_j^{A_r}(A_r/I^n_r,B_r),A)\Rightarrow
\Tor^{A_r}_{i+j}(A_r/I^n_r,A).
$$

\begin{claim}\label{cl_degenerate-at-E-two}
$\Tor_j^{A_r}(A_r/I^n_r,B_r)=0$ for every $j>0$.
\end{claim}
\begin{pfclaim} We have a change of ring spectral sequence
$$
E^2_{p,q}:=\Tor^{A_r}_p(\Tor^\Z_q(A,A_r),A_r/I_r^n)\Rightarrow
\Tor_{p+q}^\Z(A,A_r/I^n_r)
$$
and clearly $E^2_{p,q}=0$ for every $q>0$, so there follows
a natural isomorphism
$$
\Tor^{A_r}_p(B_r,A_r/I_r^n)\isom\Tor^\Z_p(A,A_r/I^n_r)
\qquad
\text{for every $p\in\N$}.
$$
However, $A_r/I^n_r$ is a free $\Z$-module, whence the claim.
\end{pfclaim}

In view of claim \ref{cl_degenerate-at-E-two} we have
$E^2_{i,j}=0$ for every $j>0$, whence a natural isomorphism
$$
\Tor_i^{B_r}(B_r/I^n_rB_r,A)\isom\Tor^{A_r}_i(A_r/I^n_r,A)
\qquad
\text{for every $i\in\N$}.
$$
Especially, $\mathrm{(a)}_\bff$ is equivalent to :
\begin{itemize}
\item[$\mathrm{(g)}_\bff$]
The inverse system $(\Tor^{B_r}_i(B_r/I_r^nB_r,A)~|~n\in\N)$
is essentially zero for every $i>0$
\end{itemize}
and $\mathrm{(a)^{un}_\bff}$ is equivalent to the corresponding
condition $\mathrm{(g)_\bff^{un}}$.

Now, set $(\bff,\bg):=(f_1,\dots,f_r,g_1,\dots,g_s)$.
In light of the foregoing we see that in order to prove
the equivalence $\mathrm{(a)_\bff\Leftrightarrow(a)_\bg}$
it suffices to check that
$$
\mathrm{(g)_\bff\Leftrightarrow(g)_{(\bff,\bg)}}
\qquad\text{and}\qquad
\mathrm{(g)_\bg\Leftrightarrow(g)_{(\bff,\bg)}}
$$
and likewise for the uniformly essentially zero variants,
so we may assume from start that $s>r$ and $f_i=g_i$ for
$i=1,\dots,r$. In this case, an easy induction on $s-r$
further reduces to the case where $s=r+1$. Now, let us say
that $g_s=\sum_{i=1}^ra_if_i$ for some $a_1,\dots,a_r\in A$,
and consider the automorphism of $A$-algebras
$\omega:B_{r+1}\isom B_{r+1}$ such that $\omega(T_i):=T_i$ for
$i=1,\dots,r$ and $\omega(T_{r+1}):=T_{r+1}-\sum_{i=1}^ra_iT_i$.
Clearly $\omega(I^n_{r+1}B_{r+1})=I^n_{r+1}B_{r+1}$ for every
$n\in\N$, and $\gamma_{\bg}\circ\omega(\bg)=
(\bff,0):=(f_1,\dots,f_r,0)$. It follows easily
that $\mathrm{(g)_\bg\Leftrightarrow(g)_{(\bff,0)}}$ and
$\mathrm{(g)^{un}_\bg\Leftrightarrow(g)^{un}_{(\bff,0)}}$, so
we may further assume that $\bg=(\bff,0)$.

Let $\pi:B_{r+1}\to B_r$ be the map of $A$-algebras such
that $\pi(T_i):=T_i$ for $i=1,\dots,r$ and $\pi(T_{r+1}):=0$,
and notice that $\gamma_{(\bff,0)}=\gamma_\bff\circ\pi$; if
we endow $B_r$ with the $B_{r+1}$-module structure given
by $\pi$, there follows, for every $n\in\N$, a change of
ring spectral sequence
$$
E(n)_{i,j}^2:=
\Tor_i^{B_r}(\Tor^{B_{r+1}}_j(B_{r+1}/I^n_{r+1}B_{r+1},B_r),A)
\Rightarrow\Tor^{B_{r+1}}_{i+j}(B_{r+1}/I^n_{r+1}B_{r+1},A)
$$
as well as a system of morphisms of spectral sequences as
in \eqref{eq_several-planes}.

\begin{claim}\label{cl_standard-calculation}
(i)\ \
$\Tor^{B_{r+1}}_j(B_{r+1}/I^{n+1}_{r+1}B_{r+1},B_r)=0$ for
every $j>1$.

(ii)\ \
The projection $B_{r+1}/I^{n+1}_{r+1}\to B_{r+1}/I^n_{r+1}$
induces the zero map
$$
\Tor^{B_{r+1}}_1(B_{r+1}/I^{n+1}_{r+1}B_{r+1},B_r)
\to\Tor^{B_{r+1}}_1(B_{r+1}/I^n_{r+1}B_{r+1},B_r)
\qquad
\text{for every $n\in\N$}.
$$
\end{claim}
\begin{pfclaim} The $B_{r+1}$-module $B_r$ admits the
free resolution $B_{r+1}\to B_{r+1}\xrightarrow{\ \pi\ }B_r$
given by scalar multiplication by $T_{r+1}$. Assertion
(i) follows immediately, and we also deduce that
$$
\Tor^{B_{r+1}}_1(B_{r+1}/I^{n+1}_{r+1}B_{r+1},B_r)\isom
\Ker(B_{r+1}/I^{n+1}_{r+1}B_{r+1}\xrightarrow{\ T_{r+1}\ }
B_{r+1}/I^{n+1}_{r+1}B_{r+1})\subset I^n_{r+1}/I^{n+1}_{r+1}.
$$
Moreover, the transition maps in (ii) are induced by
the inclusions $I^n_{r+1}\subset I^{n-1}_{r+1}$ for every
$n>0$, whence (ii).
\end{pfclaim}

From claim \ref{cl_standard-calculation} we see that
$E(n)^2_{i,j}=0$ for every $j>1$, and $(E(n)^2_{i,1}~|~n\in\N)$
is a uniformly essentially zero inverse system. Hence,
$E(n)^\infty_{i,j}$ is uniformly essentially zero for every
$j>0$, and $E(n)^\infty_{i,0}=E(n)^2_{i,0}$ for every $n,i\in\N$.
Taking into account lemma \ref{lem_inverse-Serre-subcat} we
conclude that $\mathrm{(g)_{(\bff,0)}}$ (resp.
$\mathrm{(g)^{un}_{(\bff,0)}}$) holds if and only if
$(E(n)^2_{i,0}~|~n\in\N)$ is essentially zero (resp.
uniformly essentially zero) for every $i>0$. But
$E(n)^2_{i,0}=\Tor_i^{B_r}(B_r/I^n_rB_r,A)$, so the
latter condition is precisely $\mathrm{(g)_\bff}$
(resp. $\mathrm{(g)^{un}_\bff}$).

(iii): By assumption, there exists an integer $k>0$
such that $f^k_i\in J$ and $g^k_j\in I$ for every
$i=1,\dots,r$ and every $j=1,\dots,s$. By the
foregoing proof of (iv), we deduce that
$$
\mathrm{(a)_\bff\Leftrightarrow(a)_{(\bff,\bg^k)}}
\qquad\text{and}\qquad
\mathrm{(a)_\bg\Leftrightarrow(a)_{(\bff^k,\bg)}}.
$$
Hence, we are reduced to checking that
$\mathrm{(a)_{(\bff^k,\bg)}\Leftrightarrow(a)_{(\bff,\bg)}
\Leftrightarrow(a)_{(\bff,\bg^k)}}$. We show the first
equivalence; obviously the same argument will yield
the second one. To this aim, for every integers
$n\geq m\geq 0$ let $\phi^{(n,m)}_\bullet:
\bK_\bullet(\bff^n,\bg^n)\to\bK_\bullet(\bff^m,\bg^m)$
and $\psi^{(n,m)}_\bullet:
\bK_\bullet(\bff^n,\bg^{kn})\to\bK_\bullet(\bff^m,\bg^{km})$
be the transition morphisms of the inverse systems
$(\bK_\bullet(\bff^n,\bg^n)~|~n\in\N)$ and
$(\bK_\bullet(\bff^n,\bg^{kn})~|~n\in\N)$. Then, for
every $n,k\in\N$ we have morphisms
$$
\bK_\bullet(\bff^{kn},\bg^{k^2n})\xrightarrow{\ \beta_\bullet\ }
\bK_\bullet(\bff^{kn},\bg^{kn})\xrightarrow{\ \beta'_\bullet\ }
\bK_\bullet(\bff^n,\bg^{kn})\xrightarrow{\ \beta''_\bullet\ }
\bK_\bullet(\bff^n,\bg^n)
$$
such that $\beta'_\bullet\circ\beta_\bullet=\psi^{(kn,n)}_\bullet$
and $\beta''_\bullet\circ\beta'_\bullet=\phi^{(kn,n)}_\bullet$
(see \eqref{subsec_def-bphis}), whence the contention.
\end{proof}

\begin{remark}\label{rem_logic-argument}
(i)\ \
By carefully tracking the estimates arising in the proof,
it is possible to refine proposition
\ref{prop_long-equivalences}(ii) as follows. For every
$r\in\N$ and every sequence of integers $(c(i)~|~i\in\N)$
there exists a sequence of integers $(d(i)~|~i\in\N)$, such
that the following holds : for $A$ and $\bff$ of length $r$
as in \eqref{subsec_badabum}, suppose that the inverse
system $(\Tor^{A_r}_{i+1}(A_r/I_r^n,A)~|~n\in\N)$ is uniformly
essentially zero with step $\leq c(i)$, for every $i\in\N$.
Then the inverse system $(H_{i+1}(\bff,I^n)~|~n\in\N)$ is
uniformly essentially zero of step $\leq d(i)$, for
every $i\in\N$.

(ii)\ \
A shorter alternative proof for the refinement of (i)
can be proven as follows. Fix such a sequence $c_\bullet$,
and suppose that the assertion fails; then we may find
a system of pairs $((R_\lambda,\bff_\lambda)~|~\lambda\in\N)$
consisting of a ring $R_\lambda$ and a sequence $\bff_\lambda$
of $r$ elements of $R_\lambda$, for every $\lambda$, such
that the following holds. For every $\lambda\in\N$, denote
by $\beta_\lambda:A_r\to R_\lambda$ the ring homomorphism
associated with the sequence $\bff_\lambda$, as in
\eqref{subsec_badabum}; then :
\begin{enumerate}
\alphaenu
\item
the inverse system $(\Tor^{A_r}_{i+1}(A_r/I_r^n,R_\lambda)~|~n\in\N)$
is uniformly essentially zero with step $\leq c(i)$, for
every $i,\lambda\in\N$
\item
the inverse system $(H_{i+1}(\bff,I^n)~|~n\in\N)$ is
uniformly essentially zero of step $d(i,\lambda)$, for
every $i\in\N$
\item
the sequence of integers $(d(i_0,\lambda)~|~\lambda\in\N)$
is strictly increasing, for some $i_0\in\N$.
\end{enumerate}
Now, set $J:=\oplus_{\lambda\in\N}R_\lambda$, $B:=A_r\oplus J$,
and endow $B$ with the unique ring structure such that
$(a,b)\cdot(a',b'):=(aa',ba'+a'b)$ for every $a,a'\in A_r$
and $b,b'\in J$. Then $J$ is an ideal of $B$ with $J^2=0$,
the natural projection induces a ring isomorphism
$\omega:B/J\isom A_r$, and the inclusion map $\beta:A_r\to B$
is a ring homomorphism; moreover, the $A_r$-module structure
induced via $\omega$ on $J$ agrees with the $A_r$-module
structure induced by the system of maps
$(\beta_\lambda~|~\lambda\in\N)$. We consider the sequence
$\bt:=(T_1,\dots,T_r)$ of $A_r$, and its image
$\bff:=\beta(\bt)$ in $B$, and we notice the natural
isomorphisms of $A_r$-modules
$$
\Tor^{A_r}_{i+1}(A_r/I_r^n,B)\isom
\bigoplus_{\lambda\in\N}\Tor^{A_r}_{i+1}(A_r/I_r^n,A_r)
\qquad
\text{for every $i,n\in\N$}.
$$
In light of condition (a), we deduce that $B$ satisfies
condition $\mathrm{(a)^{un}_\bff}$, so it also satisfies
$\mathrm{(c)^{un}_\bff}$, by virtue of proposition
\ref{prop_long-equivalences}(ii); on the other hand,
clearly we have as well natural isomorphisms
$$
H_i(\bff,I^n_rB)\isom H_i(\bt,I^n_r)\oplus
\bigoplus_{\lambda\in\N}H_i(\bff_\lambda,I^n_rR_\lambda)
\qquad
\text{for every $i,n\in\N$}
$$
which, in view of our assumption (c), imply that the
inverse system $(H_{i_0}(\bff,I^n_rB)~|~n\in\N)$ cannot
be uniformly essentially zero, a contradiction.
\end{remark}

\sset\subsubsection{}\label{subsec_go-to-quotient}
Let $\bff:=(f_1,\dots,f_r)$ be a finite sequence of elements
of a ring $A$, and set
$$
I:=f_1A+\cdots+f_rA
\qquad\text{and}\qquad
J:=\bigcup_{n\in\N}\Ann_A(I^n).
$$
Denote by $\bar\bff:=(\bar f_1,\dots,\bar f_r)$
the image of $\bff$ in $\bar A:=A/J$.

\begin{proposition}\label{prop_a-un-and-quots}
With the notation of \eqref{subsec_go-to-quotient},
the following conditions are equivalent :
\begin{enumerate}
\alphaenu
\item
The ring $A$ satisfies condition $\mathrm{(a)}_\bff$ (resp.
$\mathrm{(a)^{un}_\bff}$, resp. $\mathrm{(c)^{un}_{\bff}}$)
of \eqref{subsec_badabum}.
\item
$\bar A$ satisfies condition $\mathrm{(a)}_{\bar\bff}$ (resp.
$\mathrm{(a)^{un}_{\bar\bff}}$, resp $\mathrm{(c)^{un}_{\bar\bff}}$)
and there exists $k\in\N$ with $I^kJ=0$.
\end{enumerate}
\end{proposition}
\begin{proof}(a)$\Rightarrow$(b): by proposition
\ref{prop_long-equivalences}(i) the ring $A$ satisfies
condition $\mathrm{(c)}_\bff$, so there exists $k\in\N$
such that the natural map $H_r(\bff,I^k)\to H_r(\bff,A)$
is the zero morphism. In other words
$$
I^k\cap\Ann_A(I)=0.
$$
Therefore we have as well $I^k\cap J=0$, and it follows
easily that $J=\Ann_A(I^k)$. If $A$ satisfies
$\mathrm{(c)^{un}_{\bar\bff}}$, notice that for every $n\geq k$
the natural map $I^n\to I^n\bar A$ is then bijective, and
thus the same holds for the induced map
$H_i(\bff,I^n)\to H_i(\bff,I^n\bar A)$, for every $i\in\N$,
so $\bar A$ fulfills $\mathrm{(c)^{un}_{\bar\bff}}$.

If $A$ satisfies $\mathrm{(a)}_{\bar\bff}$ (resp.
$\mathrm{(a)^{un}_{\bar\bff}}$), define $A_r$ and $\beta_\bff$
as in \eqref{subsec_badabum}, and endow $A$ and $J$ with the
$A_r$-module structure induced by $\beta_\bff$; from lemma
\ref{lem_ess-zero-Tors} we deduce that the inverse system
$T(J,i)_\bullet:=(\Tor_i^{A_r}(A_r/I_r^n,J)~|~n\in\N)$ is
uniformly essentially zero for every $i>0$. By assumption,
the inverse system $T(A,i):=(\Tor_i^{A_r}(A_r/I_r^n,A)~|~n\in\N)$
is essentially zero (resp. uniformly essentially zero);
moreover, the short exact sequence of $A_r$-modules
$0\to J\to A\to\bar A\to 0$ yields exact sequences
of inverse systems :
\set\begin{equation}\label{eq_T_A_J}
T(J,i)_\bullet\to
T(A,i)_\bullet\to
T(\bar A,i)_\bullet:=(\Tor_i^{A_r}(A_r/I_r^n,\bar A)~|~n\in\N)
\to T(J,i-1)_\bullet
\end{equation}
for every $i>0$. By virtue of lemma
\ref{lem_inverse-Serre-subcat}, it follows already that
$T(\bar A,i)_\bullet$ is essentially zero (resp. uniformly
essentially zero) for every $i>1$. By the same token, the
case $i=1$ yields the exact sequence of inverse systems
$$
T(A,1)_\bullet\to T(\bar A,1)\to(J/I^nJ~|~n\in\N)
\xrightarrow{\ \phi_\bullet\ }(A/I^n~|~n\in\N).
$$
But the foregoing shows that $\phi_n$ is injective
for every $n\geq k$, hence $\Ker\,\phi_n$ is uniformly
essentially zero, and finally $T(\bar A,1)_\bullet$ is
essentially zero (resp. uniformly essentially zero),
again by lemma \ref{lem_inverse-Serre-subcat}.
Summing up, this shows that $\bar A$ fulfills
condition $\mathrm{(a)}_{\bar\bff}$ (resp.
$\mathrm{(a)^{un}_{\bar\bff}}$), as stated.

(b)$\Rightarrow$(a): Suppose first that $\bar A$ satisfies
$\mathrm{(a)}_\bff$ (resp. $\mathrm{(a)^{un}_\bff}$); by
assumption, $J$ is an $A_r/I_r^k$-module, so the inverse
system $T(J,i)_\bullet$ is uniformly essentially zero for
every $i>0$, by lemma \ref{lem_ess-zero-Tors}; the inverse
system $T(\bar A,i)_\bullet$ is essentially zero (resp.
uniformly essentially zero), and thus the same holds
also for $T(A,i)_\bullet$, in view of \eqref{eq_T_A_J}
and lemma \ref{lem_inverse-Serre-subcat}. Lastly, if
$\bar A$ satisfies $\mathrm{(c)^{un}_{\bar\bff}}$,
then it also satisfies condition $\mathrm{(a)}_\bff$, by
proposition \ref{prop_long-equivalences}(i), hence $A$
satisfies $\mathrm{(a)}_\bff$, by the foregoing case;
but as we have seen, this implies that the natural
map $H_i(\bff,I^n)\to H_i(\bff,I^n\bar A)$ is bijective
for every $i\in\N$ and sufficiently large $n$, so $A$
satisfies $\mathrm{(c)^{un}_{\bar\bff}}$.
\end{proof}

\sset\subsubsection{}\label{subsec_col-Exts}
Let $A$, $I$ and $\bff$ be as in \eqref{subsec_go-to-quotient},
and for every $n>0$ denote by $I^{(n)}\subset A$ the ideal
generated by $\bff^n:=(f_i^n~|~i=1,\dots,r)$. For all
$m\geq n>0$ we deduce natural commutative diagrams of
complexes :
$$
\xymatrix{
\bK_\bullet(\bff^m) \ar[r] \ar[d]_{\bphi_{\bff^{m-n}}} &
A/I^{(m)}[0] \ar[d]^{\pi_{mn}} \\
\bK_\bullet(\bff^n) \ar[r] & A/I^{(n)}[0]
}$$
(notation of \eqref{subsec_def-bphis}) where $\pi_{mn}$ is
the natural surjection, whence a compatible system of maps:
\set\begin{equation}\label{eq_from-Ext-to-Kos}
\Hom_{\sD(A\Mod)}(A/I^{(n)}[0],M^\bullet)\to
\Hom_{\sD(A\Mod)}(\bK_\bullet(\bff^n),M^\bullet).
\end{equation}
for every $n\geq 0$ and every complex $M^\bullet$ in
$\sD^+(A\Mod)$. Since $\bK_\bullet(\bff^n)$ is a complex
of free $A$-modules, \eqref{eq_from-Ext-to-Kos} translates
as a direct system of maps :
\set\begin{equation}\label{eq_Ext-syst}
R^i\Hom_A(A/I^{(n)},M^\bullet)\to H^i(\bff^n,M^\bullet)
\qquad\text{for all $n\in\N$ and every $i\in\Z$}.
\end{equation}
Notice that $I^{nr-r+1}\subset I^{(n)}\subset I^n$ for every
$n>0$, hence the colimit of the system \eqref{eq_Ext-syst}
is equivalent to a natural map :
\set\begin{equation}\label{eq_col-Exts}
\colim_{n\in\N}R^i\Hom^\bullet_A(A/I^n,M^\bullet)\to
\colim_{n\in\N} H^i(\bff^n,M^\bullet)
\qquad\text{for every $i\in\Z$}.
\end{equation}

\begin{lemma}\label{lem_ess-zero}
With the notation of \eqref{subsec_col-Exts}, the following
conditions are equivalent :
\begin{enumerate}
\alphaenu
\item
The map \eqref{eq_col-Exts} is an isomorphism for every
$M^\bullet\in\Ob(\sD^+(A\Mod))$ and all $i\in\Z$.
\item
$\colim_{n\in\N} H^i(\bff^n,J)=0$ for every $i>0$ and every
injective $A$-module $J$.
\item
The ring $A$ satisfies condition $\mathrm{(d)_\bff}$
of \eqref{subsec_badabum}.
\romanenu
\end{enumerate}
\end{lemma}
\begin{proof} (a)\,$\Rightarrow$\,(b) is obvious. Next, if
$J$ is an injective $A$-module, we have natural isomorphisms
\set\begin{equation}\label{eq_nat-iso-inj}
H^i(\bff^n,J)\simeq\Hom_A(H_i\bK_\bullet(\bff^n),J)
\qquad\text{for all $n\in\N$}.
\end{equation}
which easily implies that (c)\,$\Rightarrow$\,(b).

(b)\,$\Rightarrow$\,(c) : Indeed, for any $p\in\N$ let us
choose an injection $\phi:H_i\bK_\bullet(\bff^n)\to J$ into an injective
$A$-module $J$. By \eqref{eq_nat-iso-inj} we can regard $\phi$
as an element of $H^i(\bff^n,J)$; by (b) the image of
$\phi$ in $\Hom_A(H_i\bK_\bullet(\bff^q),J)$ must vanish if $q>p$
is large enough. This can happen only if
$H_i\bK_\bullet(\bff^q)\to H_i\bK_\bullet(\bff^p)$ is the zero
map.

(b)\,$\Rightarrow$\,(a) : Let $M^\bullet\to J^\bullet$ be an
injective resolution of the complex $M^\bullet$. The double
complex
$\colim_{n\in\N}\Hom^\bullet_A(\bK_\bullet(\bff^n),J^\bullet)$
determines two spectral sequences :
$$
\begin{aligned}
E_1^{pq}:=\colim_{n\in\N}\Hom_A(\bK_p(\bff^n),H^qJ^\bullet) &
\Rightarrow
\colim_{n\in\N}R^{p+q}\Hom^\bullet_A(\bK_\bullet(\bff^n),M^\bullet) \\
F_1^{pq}:=\colim_{n\in\N} H^p(\bff^n,J^q)\simeq
\colim_{n\in\N} \Hom_A(H_p\bK_\bullet(\bff^n),J^q) &
\Rightarrow\colim_{n\in\N}
R^{p+q}\Hom^\bullet_A(\bK_\bullet(\bff^n),M^\bullet).
\end{aligned}
$$
Clearly $E_1^{pq}=0$ whenever $q>0$, and (b) says that
$F_1^{pq}=0$ for $p>0$. Hence these two spectral sequences
degenerate and we deduce natural isomorphisms :
$$
\colim_{n\in\N}R^q\Hom^\bullet_A(A/I^{(n)},M^\bullet)\simeq F^{0q}_2
\isom E_2^{q0}\simeq\colim_{n\in\N}H^q(\bff^n,M^\bullet).
$$
By inspection, one sees easily that these isomorphisms are
the same as the maps \eqref{eq_col-Exts}.
\end{proof}

\begin{lemma}\label{lem_Hartsho}
In the situation of \eqref{subsec_col-Exts}, suppose that the
following holds. For every finitely presented quotient $B$ of
$A$ and every $i=1,\dots,r$, there exists $p\in\N$ such that
$$
\Ann_B(f_i^q)=\Ann_B(f_i^p)
\qquad
\text{for every $q\geq p$.}
$$
Then the ring $A$ satisfies condition $\mathrm{(d)_\bff}$
of \eqref{subsec_badabum}.
\end{lemma}
\begin{proof} We shall argue by induction on $r$. If $r=1$,
then $\bff=(f)$ for a single element $f\in A$. In this case,
our assumption ensures that there exists $p\in\N$ such that
$\Ann_A(f^q)=\Ann_A(f^p)$ for every $q\geq p$. It follows
easily that
$H_1(\bphi_{f^p}):H_1\bK_\bullet(f^{p+k})\to H_1\bK_\bullet(f^k)$
is the zero map for every $k\geq 0$ (notation of
\eqref{subsec_def-bphis}), whence the claim.

Next, suppose that $r>1$ and that the claim is known for
all sequences of less than $r$ elements. Set
$\bg:=(f_1,\dots,f_{r-1})$ and $f:=f_r$. Specializing
\eqref{eq_dist-triang} to our current situation, we derive
short exact sequences :
$$
0\to H_0(f^n,H_p\bK_\bullet(\bg^n))\to H_p\bK_\bullet(\bff^n)
\to H^0(f^n,H_{p-1}\bK_\bullet(\bg^n))\to 0
$$
for every $p>0$ and $n\geq 0$; for a fixed $p$, this is an
inverse system of exact sequences. By induction, the inverse
system $(H_i\bK_\bullet(\bg^n)~|~n\in\N)$ is essentially zero
for $i>0$; in light of lemma \ref{lem_inverse-Serre-subcat},
we deduce already that the inverse system
$(H_i\bK_\bullet(\bff^n)~|~n\in\N)$ is essentially zero
for all $i>1$. To conclude, we are thus reduced to showing
that the inverse system
$(T_n:=H^0(f^n,H_0\bK_\bullet(\bg^n))~|~n\in\N)$ is essentially
zero. However $A_n:=H_0\bK_\bullet(\bg^n)=A/(g_1^n,\dots,g^n_{r-1})$
is a finitely presented quotient of $A$ for any fixed $n\in\N$,
hence the foregoing case $r=1$ shows that the inverse system
$(T_{mn}:=\Ann_{A_n}(f^m)~|~m\in\N)$ is essentially zero.
Let $m\geq n$ be chosen so that $T_{mn}\to T_{nn}$ is the
zero map; then the composition
$T_m=T_{mm}\to T_{mn}\to T_{nn}=T_n$ is zero as well.
\end{proof}

\begin{remark}\label{rem_Hartsho}
Notice that the condition of lemma \ref{lem_Hartsho} is verified
when $A$ is noetherian. Combining with lemma \ref{lem_ess-zero},
we conclude that if $A$ is noetherian, the map \eqref{eq_col-Exts}
is an isomorphism, for every finite system $\bff$ of elements of
$A$, and every $i\in\N$.
\end{remark}

\subsection{Filtered rings and Rees algebras}
Some of the following material is borrowed from
\cite[Appendix III]{Bj}, where much more can be found.

\begin{definition}\label{def_Rees}
Let $R$ be a ring, $A$ an $R$-algebra.
\begin{enumerate}
\item
An {\em $R$-algebra filtration on $A$} is an increasing
exhaustive filtration $\Fil_\bullet A$ indexed by $\Z$ and
consisting of $R$-submodules of $A$, such that :
$$
1\in\Fil_0 A \qquad\text{and}\qquad
\Fil_i A\cdot\Fil_j A\subset\Fil_{i+j}A
\qquad\text{for every $i,j\in\Z$}.
$$
The pair $\underline A:=(A,\Fil_\bullet A)$ is called a
{\em filtered $R$-algebra}.
\item
Let $M$ be an $A$-module. An {\em $\underline A$-filtration}
on $M$ is an increasing exhaustive filtration $\Fil_\bullet M$
consisting of $R$-submodules, and such that :
$$
\Fil_iA\cdot\Fil_jM\subset\Fil_{i+j}M
\qquad\text{for every $i,j\in\Z$}.
$$
The pair $\underline M:=(M,\Fil_\bullet M)$ is called a
{\em filtered $\underline A$-module}.
\item
Let $U$ be an indeterminate. The {\em Rees algebra\/} of
$\underline A$ is the $\Z$-graded subring of $A[U,U^{-1}]$
$$
\sR(\underline A)_\bullet:=\bigoplus_{i\in\Z}U^i\cdot\Fil_i A.
$$
\item
Let $\underline M:=(M,\Fil_\bullet M)$ be a filtered
$\underline A$-module. The {\em Rees module} of $\underline M$
is the graded $\sR(\underline A)_\bullet$-module :
$$
\sR(\underline M)_\bullet:=\bigoplus_{i\in\Z}U^i\cdot\Fil_i M.
$$
\end{enumerate}
\end{definition}

\begin{lemma} Let $R$ be a ring, $\underline A:=(A,\Fil_\bullet A)$
a filtered $R$-algebra, $\gr_\bullet\underline A$ the associated
graded $R$-algebra, $\sR(\underline A)_\bullet\subset A[U,U^{-1}]$
the Rees algebra of $\underline A$. Then there are natural isomorphisms
of graded $R$-algebras :
$$
\sR(\underline A)_\bullet/U\sR(\underline A)_\bullet\simeq
\gr_\bullet\underline A
\qquad
\sR(\underline A)_\bullet[U^{-1}]\simeq A[U,U^{-1}]
$$
and of $R$-algebras :
$$
\sR(\underline A)_\bullet/(1-U)\sR(\underline A)_\bullet\simeq A.
$$
\end{lemma}
\begin{proof} The isomorphisms with $\gr_\bullet\underline A$
and with $A[U,U^{-1}]$ follow directly from the definitions.
For the third isomorphism, it suffices to remark that
$A[U,U^{-1}]/(1-U)\simeq A$.
\end{proof}

\begin{example}\label{ex_Rees}
Let $A$ be any ring, $I\subset A$ any ideal, and endow
$A$ with its $I$-adic filtration $\Fil_\bullet A$, so that
$\Fil_{-i}A:=I^i$ for every $i\in\N$ and $\Fil_iA=A$ for
every $i>0$. The pair $\underline A:=(A,\Fil_\bullet)$
is obviously a filtered $A$-algebra, and we denote by
$\sR(A,I)_\bullet$ its Rees algebra.
\end{example}

\begin{definition}\label{def_good-filtr}
Let $R$ be a ring, $\underline A:=(A,\Fil_\bullet A)$
a filtered $R$-algebra.
\begin{enumerate}
\item
Suppose that $A$ is of finite type over $R$, let $\bx:=(x_1,\dots,x_n)$
be a finite set of generators for $A$ as an $R$-algebra, and
$\bk:=(k_1,\dots,k_n)$ a sequence of $n$ integers; the
{\em good filtration $\Fil_\bullet A$ attached to the pair $(\bx,\bk)$}
is the $R$-algebra filtration such that $\Fil_iA$ is the $R$-submodule
generated by all the elements of the form
$$
\prod_{j=1}^n x_j^{a_j}\
\qquad\text{where :}\qquad
\sum_{j=1}^n a_jk_j\leq i
\quad\text{and}\quad a_1,\dots,a_n\geq 0
$$
for every $i\in\Z$. A filtration $\Fil_\bullet A$ on $A$
is said to be {\em good\/} if it is the good filtration
attached to some system of generators $\bx$ and some sequence
of integers $\bk$.
\item
The filtration $\Fil_\bullet A$ is said to be {\em positive\/}
if it is the good filtration associated with a pair $(\bx,\bk)$
as in (i), such that moreover $k_i>0$ for every $i=1,\dots,n$.
\item
Let $M$ be a finitely generated $A$-module. An
$\underline A$-filtration $\Fil_\bullet M$ is called a
{\em good filtration\/} if $\sR(M,\Fil_\bullet M)_\bullet$ is
a finitely generated $\sR(\underline A)_\bullet$-module.
\end{enumerate}
\end{definition}

\begin{example}\label{ex_Rees-free}
Let $A:=R[t_1,\dots,t_n]$ be the free $R$-algebra in $n$
indeterminates. Choose any sequence $\bk:=(k_1,\dots,k_n)$
of integers, and denote by $\Fil_\bullet A$ the good filtration
associated with $\bt:=(t_1,\dots,t_n)$ and $\bk$. Then
$\sR(A,\Fil_\bullet A)_\bullet$ is isomorphic, as a graded $R$-algebra,
to the free polynomial algebra $A[U]=R[U,t_1,\dots,t_n]$, endowed
with the grading such that $U\in\gr_1A[U]$ and $t_j\in\gr_{k_j}A[U]$
for every $j\leq n$. Indeed, a graded isomorphism can be defined
by the rule : $U\mapsto U$ and $t_j\mapsto U^{k_j}\cdot t_j$
for every $j\leq n$. The easy verification shall be left to the
reader.
\end{example}

\sset\subsubsection{}\label{subsec_good-filtr}
Let $\underline A$ and $M$ be as in definition \ref{def_good-filtr}(iii);
suppose that $\bm:=(m_1,\dots,m_n)$ is a finite system of generators
of $M$, and let $\bk:=(k_1,\dots,k_n)$ be an arbitrary sequence
of $n$ integers. To the pair $(\bm,\bk)$ we attach a filtered
$\underline A$-module $\underline M:=(M,\Fil_\bullet M)$, by
declaring that :
\set\begin{equation}\label{eq_standard-filtr}
\Fil_i M:=m_1\cdot\Fil_{i-k_1} A+\cdots+m_n\cdot\Fil_{i-k_n} A
\qquad\text{for every $i\in\Z$}.
\end{equation}
Notice that the homogeneous elements
$m_1\cdot U^{k_1},\dots,m_n\cdot U^{k_n}$ generate the graded Rees
module $\sR(\underline M)_\bullet$, hence $\Fil_\bullet M$ is
a good $\underline A$-filtration. Conversely :

\begin{lemma}\label{lem_good-filtr}
For every good filtration $\Fil_\bullet M$ on $M$, there exist
a sequence $\bm$ of generators of $M$ and a sequence of integers
$\bk$, such that $\Fil_\bullet M$ is of the form \eqref{eq_standard-filtr}.
\end{lemma}
\begin{proof} Suppose that $\Fil_\bullet M$ is a good filtration;
then $\sR(\underline M)_\bullet$ is generated by finitely many
homogeneous elements
$m_1\cdot U^{k_1},\dots,m_n\cdot U^{k_n}$. Thus,
$$
\sR(\underline M)_i:=
U^i\cdot\Fil_i M=m_1\cdot U^{k_1}\cdot A_{i-k_1}+\cdots+
m_n\cdot U^{k_n}\cdot A_{i-k_n}
\qquad\text{for every $i\in\Z$}
$$
which means that the sequences $\bm:=(m_1,\dots,m_n)$ and
$\bk:=(k_1,\dots,k_n)$ will do.
\end{proof}

\sset\subsubsection{}\label{subsec_filtered-tors}
Let $A\to B$ and $A\to C$ be maps of noetherian rings,
and suppose that either $B$ or $C$ is an $A$-algebra
of finite type; let also $M$ be a finitely generated
$B$-module, $N$ a finitely generated $C$-module, and
$I\subset B$ any ideal. Then the $A$-modules
$$
T_i:=\Tor^A_i(M,N)
$$
inherit natural $B\otimes_AC$-module structures, and
we remark :

\begin{lemma} In the situation of \eqref{subsec_filtered-tors},
the $B\otimes_AC$-module $T_i$ is finitely generated for
every $i\in\N$.
\end{lemma}
\begin{proof} Say that $B$ is an $A$-algebra of finite
type; the same argument will apply also to the case where
$C$ is an $A$-algebra of finite type. Then, we have $B=D/J$
for some free polynomial $A$-algebra $D$ of finite type,
and some ideal $J\subset D$. Hence, $M$ is also a $D$-module
of finite type; then $D$ is a noetherian ring, so we may
find a resolution $L_\bullet\to M$ consisting of free
$D$-modules of finite type, and since $D$ is a free
$A$-module, we get natural $A$-linear isomorphisms
$\omega_i:T_i\isom H_i(L_\bullet\otimes_AN)$ for every $i\in\N$.
Moreover, $L_i\otimes_AN$ is naturally a complex of
$D\otimes_AC$-modules of finite type, so $T_i$ inherits via
$\omega_i$ a structure of $D\otimes_AC$-module of finite
type. Lastly, a simple inspection shows that the latter
$D\otimes_AC$-module structure on $T_i$ agrees with the
one induced from the natural $B\otimes_AC$-module structure,
via the surjection $D\otimes_AC\to B\otimes_AC$, whence the
contention.
\end{proof}

We endow $B$ (resp. $B\otimes_AC$) with its $I$-adic
(resp. $I\cdot(B\otimes_AC)$-adic) filtration as in example
\ref{ex_Rees}, and denote by $\underline B$ (resp.
$\underline B\otimes_AC$) the resulting filtered ring.
Also, set $\Fil_nM:=\Fil_nB\cdot M$ for every $n\in\Z$;
we define a $\underline B\otimes_AC$-filtration on $T_i$,
by the rule :
$$
\Fil_nT_i:=\Img\,(\Tor^A_i(\Fil_nM,N)\to T_i)
\qquad
\text{for every $i\in\N$ and $n\in\Z$}.
$$

\begin{proposition}\label{prop_good-Tors}
The $\underline B\otimes_AC$-filtration $\Fil_\bullet T_i$
is good, for every $i\in\N$.
\end{proposition}
\begin{proof} Pick a finite set of generators $f_1,\dots,f_r$
for $I$, and consider the surjective map of graded $B$-algebras
\set\begin{equation}\label{eq_consider}
B[U,t_1,\dots,t_r]\to\sR(\underline B)_\bullet
\quad :\quad
U\mapsto U\in\sR(\underline B)_1
\quad
t_i\mapsto f_i
\qquad
\text{for $i=1,\dots,r$}
\end{equation}
(for the grading of $B[U,t_1,\dots,t_r]$ that places the
indeterminates $t_1,\dots,t_r$ in degree $-1$, and $U$ in
degree $1$). The $B\otimes_AC$-module
\set\begin{equation}\label{eq_P_i-nicoletta}
P_{i,\bullet}:=\bigoplus_{n\in\Z}\Tor^A_i(\Fil_nM,N)
\end{equation}
carries a natural structure of graded
$\sR(\underline B)_\bullet\otimes_AC$-module, whence a
graded $B\otimes_AC[U,t_1,\dots,t_r]$-module structure
as well, via \eqref{eq_consider}, and it suffices to show :

\begin{claim}\label{cl_nicoletta}
$P_{i,\bullet}$ is a finitely generated
$B\otimes_AC[U,t_1,\dots,t_r]$-module, for every $i\in\N$.
\end{claim}
\begin{pfclaim}[] Notice that
$\sR(M,\Fil_\bullet M)_\bullet$ is a finitely generated
$\sR(\underline B)_\bullet$-module; {\em a fortiori}, it
is a finitely generated graded $B[U,t_1,\dots,t_r]$-module.
By remark \ref{rem_graded-resol}(iv), we may then find
a resolution
$$
\cdots\to L_n\xrightarrow{\ d_n\ }L_{n-1}\to\cdots\to L_0
\xrightarrow{\ d_0\ }\sR(M,\Fil_\bullet M)_\bullet
$$
where each $L_n$ is a free $B[U,t_1,\dots,t_r]$-module
of finite rank, each $d_n$ is a morphism of graded
$B[U,t_1,\dots,t_r]$-modules, and the restriction of the
resolution to the summands of degree $i$ is a flat
resolution of the $B$-module $I^iM$, for every $i\in\N$.
There follows a natural isomorphism of graded
$B\otimes_AC$-modules
\set\begin{equation}\label{eq_new-B-structure}
P_{i,\bullet}\isom H_i(L_\bullet\otimes_AN)
\qquad
\text{for every $i\in\N$}
\end{equation}
and a simple inspection shows that the
$B\otimes_AC[U,t_1,\dots,t_r]$-module structure on
$P_{i,\bullet}$ deduced via \eqref{eq_new-B-structure}
agrees with the foregoing one, so the assertion follows.
\end{pfclaim}
\end{proof}

\sset\subsubsection{}\label{subsec_pro-tors}
In the situation of \eqref{subsec_filtered-tors}, set
$B_n:=B/I^n$ for every $n\in\N$. We deduce, for every
$i\in\N$, a morphism of projective systems of
$B\otimes_AC$-modules
$$
X^i_\bullet:=(B_n\otimes_B\Tor^A_i(M,N)~|~n\in\N)
\xrightarrow{\ \phi_{i,\bullet}\ }
Y^i_\bullet:=(\Tor^A_i(B_n\otimes_BM,N)~|~n\in\N)
$$
where the transition maps of $X^i_\bullet$ and $Y^i_\bullet$
are induced by the projections $B_{n+1}\to B_n$, for every
$n\in\N$, and the morphisms $\phi_{i,n}$ are induced by
the projections $M\to B_n\otimes_BM$.

\begin{corollary}\label{cor_T_i-is-good}
With the notation of \eqref{subsec_pro-tors}, we have :
\begin{enumerate}
\item
The morphism $\phi_{i,\bullet}$ is an isomorphism of
pro-$B\otimes_AC$-modules, for every $i\in\N$.
\item
The morphism $\phi_{i,\bullet}$ induces an isomorphism
of $B\otimes_AC$-modules :
$$
\lim_{n\in\N}B_n\otimes_B\Tor^A_i(M,N)\isom
\lim_{n\in\N}\Tor^A_i(B_n\otimes_BM,N)
\qquad
\text{for every $i\in\N$}.
$$
\end{enumerate}
\end{corollary}
\begin{proof}(i): The assertion means that the systems
$(\Ker\,\phi_{i,n}~|~n\in\N)$ and $(\Coker\,\phi_{i,n}~|~n\in\N)$
are essentially zero, for every $i\in\N$. We shall prove
the more precise :

\begin{claim}\label{cl_vanish-with-bounds}
For every $i\in\N$ the systems $(\Ker\,\phi_{i,n}~|~n\in\N)$
and $(\Coker\,\phi_{i,n}~|~n\in\N)$ are uniformly essentially
zero (see definition \ref{def_essentially-zero}(ii)).
\end{claim}
\begin{pfclaim} The long exact $\Tor^A_\bullet(-,N)$-sequence
arising from the short exact sequence
$$
0\to I^nM\to M\to M/I^nM\to 0
$$
yields a natural identification
$$
\Coker\,\phi_{i,n}=\Ker\,(U^n:P_{i-1,-n}\to P_{i-1,0})
$$
where $P_{i-1,\bullet}$ is defined as in \eqref{eq_P_i-nicoletta},
and $U^n$ denotes the scalar multiplication by the same element
of $\sR(\underline B)_\bullet$, for the natural
$\sR(\underline B)_\bullet$-module structure of $P_{i-1,\bullet}$.
Under this identification, the system
$(\Coker\,\phi_{i,n}~|~n\in\N)$ becomes a direct summand of the
system of $B\otimes_AC$-modules
\set\begin{equation}\label{eq_nico-again}
(\Ker\,(U^n:P_{i-1,\bullet}\to P_{i-1,\bullet})~|~n\in\N)
\end{equation}
whose transition maps are given by multiplication by $U$.
By the same token, the inverse system
$(\Ker\,\phi_{i,n}~|~n\in\N)$ is naturally identified
with the inverse system
$(\Fil_{-n}T_i/I^nT_i~|~n\in\N)$, for every $i\in\N$.
Now, it follows easily from proposition \ref{prop_good-Tors}
and lemma \ref{lem_good-filtr} that, for every $i\in\N$,
there exists $c\in\N$ such that
$$
\Fil_{-n-c}T_i\subset I^nT_i
\qquad
\text{for every $n\in\N$}
$$
whence the sought vanishing for the maps between kernels.
Lastly, since $B\otimes_AC[U,t_1,\dots,t_r]$ is noetherian,
claim \ref{cl_nicoletta} implies that there exists $c\in\N$
such that $\Ker\,U^{c+n}=\Ker\,U^c$ for every $n\in\N$.
Therefore the transition map $\Ker\,U^{c+n}\to\Ker\,U^n$
of \eqref{eq_nico-again} vanishes for every $n\in\N$,
and the proof of the claim is complete.
\end{pfclaim}

(ii) is a standard consequence of (i) : see
\cite[Prop.3.5.7]{We}.
\end{proof}

\begin{corollary}\label{cor_a-unif-for-noether}
In the situation of \eqref{subsec_badabum}, endow $A$ with
its $I$-adic filtration $\sFil_\bullet A$, and set
$B_\bullet:=\sR(A,\Fil_\bullet A)$. Denote also by $\bff_0$
the image in $B_0$ of the sequence $\bff$. Then we have :
\begin{enumerate}
\item
$B_\bullet$ satisfies condition $\mathrm{(c)}_{\bff_0}$ if
and only if $B_\bullet$ satisfies condition
$\mathrm{(c)^{un}_{\bff_0}}$, if and only if $A$ satisfies
condition $\mathrm{(c)^{un}_\bff}$.
\item
$B_\bullet$ satisfies condition $\mathrm{(a)^{un}_{\bff_0}}$
if and only if $A$ satisfies condition
$\mathrm{(a)^{un}_\bff}$.
\item
If $A$ is noetherian, then $A$ satisfies condition
$\mathrm{(a)^{un}_\bff}$.
\end{enumerate}
\end{corollary}
\begin{proof}(iii): It suffices to apply claim
\ref{cl_vanish-with-bounds} with $A$ (resp. $I$, resp. $B$,
resp. $C$) replaced by $A_r$ (resp. by $I_r$, resp. by
$A_r$, resp. by $A$).

(i): Indeed, if condition $\mathrm{(c)}_{\bff_0}$
holds, for every $i\in\N$ there exists $k\in\N$ such that
the natural map $H_i(\bff_0,I^kB_\bullet)\to H_i(\bff_0,B)$
is the zero map. However, clearly the latter is the direct
sum of the natural maps
$$
H_i(\bff,I^{n+k})\to H_i(\bff,I^n)
\qquad
\text{for every $n\in\N$}
$$
whence $\mathrm{(c)^{un}_\bff}$. By the same token, we immediately
see that $\mathrm{(c)^{un}_\bff}\Rightarrow\mathrm{(c)^{un}_{\bff_0}}$.

(ii): Define $A_r$, $I_r$, $\beta_\bff:A_r\to A$ and
$\beta_{\bff_0}:A_r\to B_\bullet$ as in \eqref{subsec_badabum}, endow
$A_r$ with its $I_r$-adic filtration $\Fil_\bullet A_r$, and set
$$
B_{r,\bullet}:=\sR(A_r,\Fil_\bullet A_r)
\qquad
C_\bullet:=B_{r,\bullet}\otimes_{A_r}A.
$$
Let also $\bt$ be the sequence $(T_1,\dots,T_r)$ of homogenous
elements of degree zero in $B_{r,\bullet}$, and $\bg$ the
image of $\bt$ under the natural map of graded rings
$B_{r,\bullet}\to C_\bullet$; since the ring $B_{r,\bullet}$ is
noetherian, (iii) tells us that it satisfies condition
$\mathrm{(a)^{un}_\bt}$, which means that we have uniformly
essentially zero inverse systems
$$
K_{r,\bullet}:=
(\Ker\,(I_r^n\otimes_{A_r}B_{r,\bullet}\to I_r^nB_{r,\bullet})~|~n\in\N)
\qquad\text{and}\qquad
T^i_{r,\bullet}:=(\Tor_i^{A_r}(I_r^n,B_{r,\bullet})~|~n\in\N)
$$
for every $i>0$. Suppose now that $\mathrm{(a)^{un}_\bff}$
holds; this means that also the inverse systems
$$
K_\bullet:=(\Ker\,(I_r^n\otimes_{A_r}A\to I^n)~|~n\in\N)
\quad\text{and}\quad
T^i_\bullet:=(\Tor_i^{A_r}(I_r^n,A)~|~n\in\N)
\quad
\text{for every $i>0$}
$$
are uniformly essentially zero; in this case, we notice :

\begin{claim}\label{cl_step-kills}
For every $i>0$ there exists $c\in\N$ such that
$I^c\cdot T^i_n=0$ for every $n\in\N$.
\end{claim}
\begin{pfclaim} By assumption, for every $i\in\N$
there exists $c\in\N$ such that the transition map
$\phi_{n+c,n}:T^i_{n+c}\to T^i_n$ is the zero morphism
for every $n\in\N$. Now, if $a\in I^c$, the scalar
multiplication by $a$ on $I^n$ factors as the composition
of an $A$-linear map $I^n\to I^{n+c}$ and the inclusion map
$I^{n+c}\to I^n$. It follows that scalar multiplication by
$a$ on $T^i_n$ factors as the composition of an $A$-linear
map $T^i_n\to T^i_{n+c}$ and $\phi_{n+c,n}$, whence the claim.
\end{pfclaim}

\begin{claim}\label{cl_reduce-to-C}
If $C_\bullet$ satisfies condition $\mathrm{(a)^{un}_\bg}$,
then $B_\bullet$ satisfies condition $\mathrm{(a)^{un}_{\bff_0}}$.
\end{claim}
\begin{pfclaim} We have a natural surjective homomorphism
$f_\bullet:C_\bullet\to B_\bullet$ of graded rings, such that
$\Ker\,f_\bullet=\oplus_{n\in\N}K_n$, and arguing as in the
proof of claim \ref{cl_step-kills}, we see that there
exists $c\in\N$ such that $I^c\cdot K_n=0$ for every
$n\in\N$. Set
$$
J_C:=\bigcup_{n\in\N}\Ann_{C_\bullet}(I^n)
\qquad
J_B:=\bigcup_{n\in\N}\Ann_{B_\bullet}(I^n)
$$
and $\bar C:=C_\bullet/J_C$; it follows that the projection
$C_\bullet\to\bar C_\bullet$ factors through a surjective ring
homomorphism $B_\bullet\to\bar C$ whose kernel equals $J_B$,
so that $f_\bullet(J_C)=J_B$. Suppose now that $C_\bullet$
satisfies condition $\mathrm{(a)^{un}_\bg}$, and let
$\bar\bg$ be the image of $\bg$ in $\bar C$; then there
exists $k\in\N$ such that $I^kJ_C=0$, and $\bar C$
satisfies condition $\mathrm{(a)^{un}_{\bar\bg}}$
(proposition \ref{prop_a-un-and-quots}). Hence $I^kJ_B=0$,
and the assertion follows, again by proposition
\ref{prop_a-un-and-quots}.
\end{pfclaim}

In view of claim \ref{cl_reduce-to-C}, we are reduced
to checking that for every $i>0$ the inverse systems
$$
K'_\bullet:=
(\Ker\,(I_r^n\otimes_{A_r}C_\bullet\to I_r^nC_\bullet)~|~n\in\N)
\qquad\text{and}\qquad
T'^i_\bullet:=(\Tor_i^{A_r}(I_r^n,C_\bullet)~|~n\in\N)
$$
are uniformly essentially zero. However, notice that
the natural map
$$
K_{r,\bullet}\otimes_{A_r}A\to K'_\bullet
$$
is an epimorphism of inverse systems, so we get already
the assertion for $K'_\bullet$. Next, we remark:

\begin{claim}\label{cl_triple-der-prod}
The inverse system
$(H_i(A\derotimes_{A_r}B_{r,\bullet}\derotimes_{A_r}I^n_r)~|~n\in\N)$
is uniformly essentially zero for every $i>0$.
\end{claim}
\begin{pfclaim} We have a spectral sequence :
$$
E(n)^2_{pq}:=\Tor^{A_r}_p(A,\Tor^{A_r}_q(B_{r,\bullet},I^n_r))
\Rightarrow
H_{p+q}(A\derotimes_{A_r}B_{r,\bullet}\derotimes_{A_r}I^n_r)
\qquad
\text{for every $n\in\N$}
$$
hence it suffices to check that the inverse system
$(E(n)^2_{pq}~|~n\in\N)$ is uniformly essentially zero
for every $p,q\in\N$ such that $p+q>0$. The latter
assertion is already known for every $q>0$, since
the system $T^q_{r,\bullet}$ is uniformly essentially
zero. For $q=0$, we consider the short exact sequence
of inverse systems
$$
0\to K_{r,\bullet}\to(B_{r,\bullet}\otimes_{A_r}I^n_r~|~n\in\N)
\to(I^n_rB_{r,\bullet}~|~n\in\N)\to 0
$$
whence, for every $p\in\N$ an exact sequence of inverse systems
$$
X^p_\bullet:=
\Tor_p^{A_r}(A,K_{r,n}~|~n\in\N)\to(E(n)_2^{p0}~|~n\in\N)
\to Y^p_\bullet:=(\Tor_p^{A_r}(A,I^n_rB_{r,\bullet})~|~n\in\N).
$$
Now, the system $X^p_\bullet$ is uniformly essentially zero
for every $p\in\N$, since the same holds for $K_{r,\bullet}$;
moreover, $Y^p_\bullet$ is uniformly essentially zero for
$p>0$, since the same holds for $T^p_\bullet$. Hence,
$(E(n)_2^{p0}~|~n\in\N)$ is uniformly essentially zero
for $p>0$, and the proof is concluded.
\end{pfclaim}

Lastly, for every $n\in\N$ let us consider the spectral
sequence
$$
E(n)_{pq}^2:=\Tor_p^{A_r}(\Tor^{A_r}_q(A,B_{r,\bullet}),I^n_r)
\Rightarrow
H_{p+q}(A\derotimes_{A_r}B_{r,\bullet}\derotimes_{A_r}I^n_r)
$$
and notice that $(E(n)^2_{p0}~|~n\in\N)=T'^p_\bullet$ for
every $p\in\N$. To conclude, it suffices therefore to show :

\begin{claim} The inverse system $(E(n)^k_{p0}~|~n\in\N)$
is uniformly essentially zero, for every $k\geq 2$ and
every $p>0$.
\end{claim}
\begin{pfclaim} We fix $p>0$ and argue by descending induction
on $k$. Notice that $E(n)^{p+1}_{p0}=E(n)^\infty_{p0}$ for every
$n\in\N$, so claim \ref{cl_triple-der-prod} implies that the
assertion holds whenever $k\geq p+1$. Next, let $2\leq i\leq p$
and suppose that the assertion is known for every $k>i$; we
have an exact sequence of inverse systems
$$
(E(n)^{i+1}_{p0}~|~n\in\N)\to(E(n)^i_{p0}~|~n\in\N)\to
(E(n)^i_{p-i,i-1}~|~n\in\N).
$$
In view of our inductive assumption (and of lemma
\ref{lem_inverse-Serre-subcat}), we are therefore reduced
to showing that the inverse system $(E(n)^i_{p-i,i-1}~|~n\in\N)$
is uniformly essentially zero. However, by claim
\ref{cl_step-kills} there exists $c\in\N$ such that
$I^c_r\cdot\Tor^{A_r}_{i-1}(A,B_{r,\bullet})=0$, and then
the contention follows from lemma \ref{lem_ess-zero-Tors}.
\end{pfclaim}

Conversely, if $\mathrm{(a)^{un}_{\bar\bff_0}}$ holds for
$B_\bullet$, then we get immediately condition
$\mathrm{(a)^{un}_{\bar\bff}}$ for $A$, since the latter
is a direct summand of the $A_r$-module $B_\bullet$.
\end{proof}

\begin{remark}\label{rem_a-unif-for-noether}
By carefully tracking the estimates in the proof, it
is easily seen that for every sequence of integers
$(c(i)~|~i\in\N)$ there exists a sequence of integers
$(d(i)~|~i\in\N)$ such that the following refinement
holds. Suppose that for every $i\in\N$ the inverse
system $(\Tor^{A_r}_i(A_r/I^n_r,A)~|~n\in\N)$ is uniformly
essentially zero with step $\leq c(i)$; then the inverse
system $(\Tor^{A_r}_i(A_r/I^n_r,B_\bullet)~|~n\in\N)$ is
uniformly essentially zero with step $\leq d(i)$.
\end{remark}

\sset\subsubsection{}\label{subsec_lift-compltely-sec}
Let $(A,\Fil_\bullet A)$ be a filtered ring, {\em i.e.}
a filtered $\Z$-algebra, in the sense of definition
\ref{def_Rees}(i), and denote by $\gr_\bullet A$ the
associated graded ring. Suppose that the filtration
$\Fil_\bullet A$ is {\em exhaustive} and {\em separated},
that is
$$
\bigcup_{n\in\Z}\Fil_nA=A
\qquad\text{and}\qquad
\bigcap_{n\in\Z}\Fil_nA=0.
$$
The filtration $\Fil_\bullet A$ defines a linear topology
$\cT$ on the $\Z$-module $A$, and we suppose moreover
that $\cT$ is complete (and notice that $\cT$ is obviously
separated). Furthermore, let $\bff:=(f_1,\dots,f_r)$
be a finite sequence of elements of $A$; for every
$i=1,\dots,r$, pick $n_i\in\Z$ such that $f_i\in\Fil_{n_i}A$,
denote by $\bar f_i\in\gr_{n_i}A$ the image of $f_i$, and let
$\bar\bff:=(\bar f_1,\dots,\bar f_r)$ be the resulting
sequence of elements of $\gr_\bullet A$.

\begin{proposition}\label{prop_lift-compl-secant}
In the situation of \eqref{subsec_lift-compltely-sec},
suppose that the sequence $\bar\bff$ is completely secant
(see definition {\em\ref{def_compl-sec}}). Then
the same holds for the sequence $\bff$.
\end{proposition}
\begin{proof} We define as follows a filtration on the
Koszul complex $\bK_\bullet(\bff)=\bK_\bullet(\bff,A)$, {\em i.e.}
a compatible system of filtrations on the exterior powers
$\Lambda^i_A(A^{\oplus r})$ for $i=1,\dots,r$. To this aim,
for every $i\leq r$ denote by $S_i$ the set of all strictly
increasing maps $\{1,\dots,i\}\to\{1,\dots,r\}$; for every
$\phi\in S_i$ set
$$
e_\phi:=e_{\phi(1)}\wedge\cdots\wedge e_{\phi(i)}
\qquad\text{and}\qquad
n_\phi:=n_{\phi(1)}+\cdots+n_{\phi(i)}.
$$
We let $\Fil_\bullet\Lambda^i_A(A^{\oplus r})$ be the good
filtration associated, as in \eqref{subsec_good-filtr},
with the basis $(e_\phi~|~\phi\in S_i)$ and the sequence
of integers $(n_\phi~|~\phi\in S_i)$. The explicit
description in \eqref{sec_koszul-alg} shows that
the differentials of the Koszul complex are morphisms
of filtered abelian groups, for the resulting filtrations
on the terms $\bK_i(\bff)$. Thus, we have the sought
filtration $\Fil_\bullet\bK_\bullet(\bff,A)$, and we
let $\gr_\bullet\bK_\bullet(\bff,A)$ be the associated
complex of graded abelian groups; with this notation, a
simple inspection shows that we have a natural isomorphism
of complexes of abelian groups
\set\begin{equation}\label{eq_from-f-to-fbar}
\gr_\bullet\bK_\bullet(\bff,A)\isom
\bK_\bullet(\bar\bff,\gr_\bullet A).
\end{equation}
For every $j,k\in\Z$ with $j\geq k$ we set
$$
\bK^{(j,k)}_\bullet(\bff,A):=
\Fil_j\bK_\bullet(\bff,A)/\Fil_{k-1}\bK_\bullet(\bff,A).
$$
We endow $\bK^{(j,k)}_\bullet(\bff,A)$ with the
filtration induced by $\Fil_\bullet\bK_\bullet(\bff,A)$, and
we let $\gr_\bullet\bK^{(j,k)}_\bullet(\bff,A)$ be the complex
of graded abelian groups associated with this filtered
subquotient of $\bK_\bullet(\bff,A)$. According to remark
\ref{rem_low-terms}(i), we have a $1$-spectral sequence
$$
E^{pq}_1:=H^{p+q}(\gr_{-p}\bK^{(j,k)}_\bullet(\bff,A))
\Rightarrow H^{p+q}\bK^{(j,k)}_\bullet(\bff,A).
$$
However, \eqref{eq_from-f-to-fbar} and the assumption
on $\bar\bff$ imply that $E^{pq}_1=0$ whenever $p+q<0$,
and then proposition \ref{prop_convergence-simple}
implies that $H^n(\bK^{(j,k)}_\bullet(\bff,A))=0$ for
every $j,k,n$ with $j\geq k$ and $n<0$. Set as well
$\bK^{(k)}_\bullet(\bff,A):=
\bK_\bullet(\bff,A)/\Fil_{k-1}\bK_\bullet(\bff,A)$ for
every $k\in\Z$; it follows easily that
$$
H^n\bK^{(k)}_\bullet(\bff,A)=0
\qquad
\text{for every $n<0$ and every $k\in\Z$}.
$$
Next, since $\bK_n(\bff,A)$ is a free $A$-module of
finite rank for every $n\in\Z$, and since the filtration
$\Fil_\bullet A$ is separated and defines a complete
topology, we get a natural isomorphism of abelian groups
$$
\bK_n(\bff,A)\isom
\lim_{k\in\Z}\bK^{(k)}_n(\bff,A)
\qquad
\text{for every $n\in\Z$}.
$$
Lastly, in view of \cite[Prop.3.5.7, Th.3.5.8]{We} we
conclude that $H^n\bK_\bullet(\bff,A)=0$ for every
$n<0$, {\em i.e.} the sequence $\bff$ is completely
secant, as stated.
\end{proof}

\subsection{Some homotopical algebra}\label{sec_homotopy}
The methods of homotopical algebra allow to construct
derived functors of non-additive functors, in a variety
of situations. In this section, we present the basics
of this theory, beginning with its cornerstone, the
{\em standard resolution associated with a triple}, as
explained in the following paragraph.

\sset\subsubsection{}\label{subsec_cotriples}
A {\em triple} $(\top,\eta,\mu)$ on a category $\cC$ is the datum
of a functor $\top:\cC\to\cC$ together with natural transformations :
$$
\etab:\one_\cC\Rightarrow\top
\qquad
\bmu:\top\circ\top\Rightarrow\top
$$
such that the following diagrams commute :
$$
\xymatrix{
(\top\circ\top)\circ\top \ar@{=>}[d]_{\bmu*\top} \ar@{=}[r] &
\top\circ(\top\circ\top) \ar@{=>}[r]^-{\top*\bmu}  &
\top\circ\top \ar@{=>}[d]^\bmu &
\top \ar@{=>}[rr]^-{\top*\etab} \ar@{=}[drr]_{\one_\top} & &
\top\circ\top \ar@{=>}[d]^\bmu & &
\top \ar@{=>}[ll]_-{\etab*\top} \ar@{=}[dll]^{\one_\top} \\
\top\circ\top \ar@{=>}[rr]^-\bmu & & \top
& & & \top.
}$$
Dually, a cotriple $(\perp,\beps,\bdelta)$ in a category $\cC$ is
a functor $\perp:\cC\to\cC$ together with natural transformations :
$$
\beps:\perp\Rightarrow\one_\cC
\qquad
\bdelta:\perp\Rightarrow\perp\circ\perp
$$
such that the following diagrams commute :
$$
\xymatrix{ \perp \ar@{=>}[rr]^-\bdelta \ar@{=>}[d]_\bdelta & &
\perp\circ\perp \ar@{=>}[d]^{\bdelta*\perp}
& & & \perp \ar@{=}[drr]^{\one_\perp}
      \ar@{=}[dll]_{\one_\perp} \ar@{=>}[d]^\bdelta \\
\perp\circ\perp \ar@{=>}[r]^-{\perp*\bdelta} &
\perp\circ(\perp\circ\perp) \ar@{=}[r] & (\perp\circ\perp)\circ\perp
& \perp & &
\perp\circ\perp \ar@{=>}[ll]_-{\perp*\beps} \ar@{=>}[rr]^-{\beps*\perp}
& & \perp.
}$$
Notice that a cotriple in $\cC$ is the same as a triple in $\cC^o$.

\sset\subsubsection{}\label{subsec_simpl-cotriple}
A cotriple $\perp$ on $\cC$ and an object $A$ of $\cC$ determine a
simplicial object $\perp\!\!A[\bullet]$ in $\cC$; namely, for every
$n\in\N$ set $\perp A[n]:=\perp^{n+1}A$, and define face and
degeneracy operators :
$$
\begin{aligned}
\partial_i:=(\perp^i*\beps*\perp^{n-i})_A \quad &: \quad
\perp A[n]\to\perp A[n-1] \\
\sigma_i:=(\perp^i*\bdelta*\perp^{n-i})_A \quad &: \quad
\perp A[n]\to\perp A[n+1].
\end{aligned}
$$
Using the foregoing commutative diagrams, one verifies easily
that the simplicial identities \eqref{eq_simpl-identities} hold.
Moreover, the morphism $\beps_A:\perp\!A\to A$ defines
an augmentation $\perp\!\!A[\bullet]\to A$.

Dually, a triple $\top$ and an object $A$ of $\cC$ determine a
cosimplicial object $\top A[\bullet]:=\top^{\bullet+1}A$, such that
$$
\begin{aligned}
\partial^i:=(\top^i*\etab*\top^{n-i})_A \quad &: \quad \top^nA\to\top^{n+1}A \\
\sigma^i:=(\top^i*\bmu*\top^{n-i})_A \quad &: \quad \top^{n+2}A\to\top^{n+1}A
\end{aligned}
$$
which is augmented by the morphism $\etab_A:A\to\top A$.
Clearly, the rule : $A\mapsto\perp\!\!A[\bullet]$ (resp.
$A\mapsto\top A[\bullet]$) defines a functor
$$
\cC\to\hat s.\cC
\qquad
\text{(resp. $\cC\to\hat c.\cC$)}.
$$

\sset\subsubsection{}\label{subsec_perp-adj}
An adjoint pair $(G,F)$ as in \eqref{subsec_adj-pair}, with
its unit $\etab$ and counit $\beps$, determines a triple
$(\top,\etab,\bmu)$, where :
$$
\top:=F\circ G:\cB\to\cB
\qquad
\bmu:=F*\beps*G:\top\circ\top\Rightarrow\top
$$
as well as a cotriple $(\perp,\beps,\bdelta)$, where :
$$
\perp\: :=G\circ F:\cA\to\cA \qquad
\bdelta:=G*\etab*F:\:\perp\:\Rightarrow\:\perp\circ\perp.
$$
Indeed, the naturality of $\beps$ and $\etab$ easily implies
the commutativity of the diagrams of \eqref{subsec_cotriples}.
The following proposition explains how to use these triples
and cotriples to construct canonical resolutions.

\begin{proposition}\label{prop_triple-res}
In the situation of \eqref{subsec_perp-adj}, the following holds :
\begin{enumerate}
\item
For every $A\in\Ob(\cA)$ and $B\in\Ob(\cB)$, the augmented
simplicial objects :
$$
\perp\!GB[\bullet]\xrightarrow{\ \beps_{GB}\ }GB
\qquad
F\!\perp\!\!A[\bullet]\xrightarrow{\ F*\beps_A\ }FA
$$
are homotopically trivial (notation of \eqref{subsec_simpl-cotriple}).
\item
Dually, the same holds for the augmented cosimplicial objects :
$$
GB\xrightarrow{\ G*\etab_B\ }G\top B[\bullet]
\qquad
FA\xrightarrow{\ \etab_{FA}\ }\top FA[\bullet].
$$
\end{enumerate}
\end{proposition}
\begin{proof} Notice that
$$
F\!\perp\!\!A[n]=\top FA[n]
\qquad
\text{for every $[n]\in\Ob(\Delta^\wedge)$}.
$$
Therefore, for every morphism $\phi:[n]\to[m]$ in $\Delta^\wedge$,
we have two morphisms
$$
F\!\perp\!\!A[\phi]:F\!\perp\!\!A[m]\to F\!\perp\!\!A[n]
\qquad
\top FA[\phi]:F\!\perp\!\!A[n]\to F\!\perp\!\!A[m].
$$
Now, recall that the augmentation $\eps_A$ defines a natural
morphism
$$
\perp\!\!A[\eps^{\bullet+1}_{-1,0}]:\perp\!\!A[\bullet]\to s.A
\qquad
[n]\mapsto\perp\!\!A[\eps^{n+1}_{-1,0}]
$$
(notation of remark \ref{rem_path-spaces}(iv)), and one sees
that the system $(\top FA[\eps^{n+1}_{-1,0}]~|~[n]\in\Ob(\Delta))$
defines a morphism
$$
\top FA[\eps^{\bullet+1}_{-1,0}]:s.FA\to F\!\perp\!\!A[\bullet]
$$
which is right inverse to $F\!\perp\!\!A[\eps^{\bullet+1}_{-1,0}]$.
For every $n,k\in\N$ such that $k\leq n+1$, set
$$
u_{n,k}:=
(\top FA[\eps^k_{n-k,0}])\circ(F\!\perp\!\!A[\eps^k_{n-k,0}])
$$
(notation of example \ref{ex_face-and-deg}(ii)). 

\begin{claim}\label{cl_Deligne-homotopy}
The system $u_{\bullet\bullet}$ defines a homotopy from
$\one_{F\!\perp\!A[\bullet]}$ to
$(\top FA[\eps^{\bullet+1}_{-1,0}])\circ(F\!\perp\!\!A[\eps^{\bullet+1}_{-1,0}])$
(see \eqref{subsec_homotopy-by-us}).
\end{claim}
\begin{pfclaim} By unwinding the definitions, we see that
$F\!\perp\!\!A[\eps^k_{n-k,0}]$ is the composition
$$
F(\beps_{\perp^{n+1-k}A}\circ\cdots\circ\beps_{\perp^{n-1}A}\circ\beps_{\perp^n A}):
F\!\perp^{n+1}\!\!A\to F\!\perp^{n+1-k}\!\!A
\qquad
\text{for $k>0$}
$$
and $\top FA[\eps^k_{n-k,0}]$ is the composition
$$
\etab_{F\perp^nA}\circ\etab_{F\perp^{n-1}A}\circ\cdots\circ\etab_{F\perp^{n+1-k}A}:
F\!\perp^{n+1-k}\!\!A\to F\!\perp^{n+1}\!\!A
\qquad
\text{for $k>0$}
$$
and both equal $\one_{F\perp^{n+1}A}$ for $k=0$. Now, since
the face operator $\partial_i$ of $F\!\perp\!\!A[n]$ is
$F\!\perp^i(\beps_{\perp^{n-i}A})$ for every $i\leq n$, the
naturality of $\etab$ yields a commutative diagram
$$
{\diagram
F\!\perp^n\!\!A \ar[rr]^-{\partial_{i-1}} \ar[d]_{\etab_{F\!\perp^n\!A}} & &
F\!\perp^{n-1}\!\!A \ar[d]^{\etab_{F\!\perp^{n-1}\!A}} \\
F\!\perp^{n+1}\!\!A \ar[rr]^{\partial_i} & & F\!\perp^n\!\!A
\enddiagram}
\qquad
\text{for every $i>0$ and every $n\in\N$}
$$
whereas $\partial_0\circ\etab_{F\!\perp^n\!A}=\one_{F\!\perp^n\!A}$
for every $n\in\N$; whence the identities :
$$
\partial_i\circ\top FA[\eps^k_{n-k,0}]=
\left\{\begin{array}{ll}
        \top FA[\eps^{k-1}_{n-k,0}] & \quad\text{for $i<k$} \\
        \top FA[\eps^k_{n-k-1,0}]\circ\partial_{i-k} & \quad\text{for $i\geq k$}.
       \end{array}\right.
$$
On the other hand, we have the simplicial identities :
$$
\eps_i\circ\eps^k_{n-k,0}=
                   \left\{\begin{array}{ll}
                    \eps^{k+1}_{n-k} & \quad\text{for $i<k$} \\
                    \eps^k_{n+1-k}\circ\eps_{i-k} & \quad\text{for $i\geq k$}.
                          \end{array}\right.
$$
Summing up, we conclude that
$$
\partial_i\circ u_{n,k}=
                   \left\{\begin{array}{ll}
                       u_{n-1,k-1}\circ\partial_i & \quad\text{for $i<k$}\\
                       u_{n-1,k}\circ\partial_i & \quad\text{for $i\geq k$}.
                          \end{array}\right.
$$
Arguing likewise, one deduces as well the corresponding commutation
rule -- spelled out in \eqref{subsec_homotopy-by-us} -- for the
degeneracies $\sigma_i$, and the claim follows.
\end{pfclaim}

Likewise, notice that
$$
G\top B[n]=\perp\!\!GB[n]
\qquad
\text{for every $[n]\in\Delta^\wedge$}.
$$
Therefore, for every morphism $\phi:[n]\to[m]$ in $\Delta^\wedge$,
we have as well two morphisms
$$
\perp\!\!GB[\phi]:\perp\!\!GB[m]\to\perp\!\!GB[n]
\qquad
G\top B[\phi]:\perp\!\!GB[n]\to\perp\!\!GB[m]
$$
and the system $(G\top B[\eps^{n+1}_{-1,0}]~|~[n]\in\Ob(\Delta^\wedge))$
defines a morphism
$$
G\top B[\eps^{\bullet+1}_{-1,0}]:s.GB\to\perp\!\!GB[\bullet]
$$
which is right inverse to the morphism $\perp\!\!GB[\eps^{\bullet+1}_{-1,0}]$
deduced from the augmentation. Arguing as in the proof of
claim \ref{cl_Deligne-homotopy}, one checks that the rule
$$
v_{n,k}:=(G\top B[\eps^{k\vee}_{n-k,0}])\circ(\perp\!\!GB[\eps^{k\vee}_{n-k,0}])
\qquad
\text{for every $n,k\in\N$ such that $k\leq n+1$}
$$
yields a homotopy $v$ from $\one_{\perp GB[\bullet]}$ to
$(G\top B[\eps^{\bullet+1}_{-1,0}])\circ(\perp\!\!GB[\eps^{\bullet+1}_{-1,0}])$.

The dual statements admit the dual proof.
\end{proof}

\sset\subsubsection{}
Now, denote by $\Set_\circ$ the category of {\em pointed sets},
whose objects are all the pairs $(S,s)$ where $S$ is a small
set, and $s\in S$ is any element; the morphisms
$f:(S,s)\to(S',s')$ are the mappings $f:S\to S'$ such that
$f(s)=s'$. Next, let $A$ be any ring, and consider the
{\em pointed forgetful functor}
$$
\sF:A\Mod\to\Set_\circ
\qquad
M\mapsto(M,0_M)
$$
(where $0_M\in M$ is the zero element of $M$). The functor
$\sF$ admits the left adjoint
$$
\sL:\Set_\circ\to A\Mod
\qquad
(S,s)\mapsto A^{(S)}/As
$$
(where $A^{(S)}$ denotes the free $A$-module with basis
given by $S$, so $As\subset A^{(S)}$ is the direct summand
generated by the basis element $s\in S$). According to
proposition \ref{prop_triple-res}, the adjoint pair
$(\sL,\sF)$ yields a functor
$$
\perp^A_\bullet:A\Mod\to\hat s.A\Mod
$$
into the category of augmented simplicial $A$-modules
(notation of \eqref{subsec_augment-obj}), such that
$$
\sF\perp^A_\bullet\!M\to\sF M
$$
is a homotopically trivial augmented pointed simplicial set, for
every $A$-module $M$. This functor is obtained by iterating the
functor $\perp^A:=\sL\circ\sF$, so in each degree it consists
of a free $A$-module. We call $\perp^A_\bullet\!M$ the
{\em standard free resolution of $M$}.

\begin{remark}\label{rem_trivial-to-trival}
Notice that $\perp^A_\bullet\!0=s.0$, the constant simplicial
$A$-module associated with the trivial $A$-module (see
\eqref{subsec_simplicial-object}). This is the reason why we
prefer the adjoint pair $(\sL,\sF)$, rather than the analogous
pair considered in \cite[I.1.5.5]{Il}, arising from the
forgetful functor from $A$-modules to non-pointed sets.
\end{remark}

\sset\subsubsection{}\label{subsec_derive-non-add}
Let $R$ be a simplicial $A$-algebra, and define the category
$R\Mod$ of $R$-modules, as in \cite[\S8.1]{Ga-Ra}. Notice
that any morphism $S\to R$ of simplicial rings induces a
base change functor
$$
S\Mod\to R\Mod
\qquad
(M[n]~|~n\in\N)\mapsto(R[n]\otimes_{S[n]}M[n]~|~n\in\N)
$$
which is left adjoint to the forgetful functor (details left
to the reader).

Now, by applying the functors $\perp^{R[n]}_\bullet$ to the
terms $M[n]$ of a given $R$-module $M$, we obtain an augmented
simplicial $R$-module
\set\begin{equation}\label{eq_augmentation-perp}
\perp^R_\bullet\!M\to M
\end{equation}
which amounts to a bisimplicial $A$-module, whose columns
$\perp^{R[n]}_\bullet\!M[n]\to M[n]$ are augmented simplicial
$R[n]$-modules, for every $n\in\N$. Likewise, the row
$\perp^R_n\!M$ is an $R$-module, for every $n\in\N$.

\begin{lemma}\label{lem_Whitegead}
In the situation of \eqref{subsec_derive-non-add}, let
$\phi:M\to N$ be any quasi-isomorphism of $R$-modules.
Then the induced morphism
$$
\perp^R_n\!\phi:\perp^R_n\!M\to\perp^R_n\!N
$$
is a homotopy equivalence, for every $n\in\N$.
\end{lemma}
\begin{proof} Set $\perp^R_{-1}\!M:=M$, and notice the natural
isomorphisms
\set\begin{equation}\label{eq_induce-from-Z}
\perp_n^R\!M\isom R\otimes_{s.\Z}\perp^{s.\Z}\circ\perp_{n-1}^R\!M
\qquad
\text{for every $n\in\N$}
\end{equation}
(where $s.\Z$ is the constant simplicial ring arising from
$\Z$ (see \eqref{subsec_simplicial-object}), which is the
initial object in the category of simplicial rings). The
assumption means that $\perp^R_{-1}\phi$ is a quasi-isomorphism.
However, \eqref{eq_induce-from-Z} and Whitehead's theorem
(\cite[I.2.2.3]{Il}) imply that if -- for a given $n\in\N$ --
the map $\perp^R_{n-1}\phi$ is a quasi-isomorphism, then
$\perp^R_n\phi$ is a homotopy equivalence. We may thus
conclude by a simple induction on $n$.
\end{proof}

\begin{example}\label{ex_der-tens-products}
Let $M$ and $N$ two $R$-modules. We may define two functors
$$
M\ellotimes_R-
\quad
\text{(\ resp. $-\ellotimes_RN$\ )}
\quad : \quad
R\Mod\to R\Mod
$$
by the rules :
$$
P\mapsto M\otimes_R(\perp_\bullet^R\!P)^\Delta
\quad
\text{(\ resp. $P\mapsto(\perp_\bullet^R\!P)^\Delta\otimes_RN$\ )}
\qquad
\text{for every $R$-module $P$}
$$
where $\Delta$ is the diagonal functor, from $\cC$-bisimplicial
to $\cC$-simplicial objects (see \eqref{subsec_bisimplex}). 

(i)\ \
We claim that these functors descend to the derived category.
That is, suppose that $\phi:P\to P'$ is a quasi-isomorphism;
then $M\ellotimes_R\phi$ and $\phi\ellotimes_RN$ are both
quasi-isomorphisms. Indeed, lemma \ref{lem_Whitegead} implies
that the induced morphisms
$$
M\otimes_R\perp_n^R\!\phi:
M\otimes_R\perp_n^R\!P\to M\otimes_R\perp_n^R\!P'
$$
are quasi-isomorphisms, for every $n\in\N$, and likewise
for $\phi\otimes_R\perp_n^R\!N$, so the assertion follows
from Eilenberg-Zilber's theorem \ref{th_Eilenberg-Zilber}(i).

(ii)\ \
Also, we claim that the notation $M\ellotimes_RN$ is
unambiguous. That is, we may compute this object by
applying the functor $M\ellotimes_R-$ to $N$, or by
applying the functor $-\ellotimes_RN$ to $M$, and the
resulting two $R$-modules are naturally isomorphic in
$\sD(R\Mod)$. Indeed, it suffices to check that the
natural morphisms
$$
M\otimes_R(\perp^R_\bullet\!N)^\Delta\xleftarrow{\ \alpha\ }
(\perp^R_\bullet\!M)^\Delta\otimes_R(\perp^R_\bullet\!N)^\Delta
\xrightarrow{\ \beta\ }(\perp^R_\bullet\!M)^\Delta\otimes_RN
$$
induced by the augmentation \eqref{eq_augmentation-perp},
are both quasi-isomorphisms. We check this for $\alpha$;
the same argument shall apply to $\beta$. To ease notation,
set $L:=(\perp^R_\bullet\!N)^\Delta$. Then $\alpha$ is
deduced by extracting the diagonal from the augmented
simplicial $R$-module
\set\begin{equation}\label{eq_check-column}
\perp^R_\bullet\!M\otimes_RL\to M\otimes_RL
\end{equation}
(so, the $n$-th column of $s.L$ is a constant and flat
simplicial $R[n]$-module, for every $n\in\N$). Thus,
we are reduced to checking that the columns of
\eqref{eq_check-column} are aspherical. However,
for every $n\in\N$, the $n$-th column is the augmented
simplicial $R[n]$-module
$(\perp^{R[n]}_\bullet\!M[n])\otimes_{R[n]}L[n]$; since
$L[n]$ is flat, it then suffices to recall that the
standard free resolution is aspherical.

(iii)\ \
The foregoing discussion yields a well defined functor
$$
-\ellotimes_R-:\sD(R\Mod)\times\sD(R\Mod)\to\sD(R\Mod)
$$
called the {\em left derived tensor product}. The same method
shall be applied hereafter to construct derived functors
of certain non-additive functors.
\end{example}

\begin{remark}\label{rem_der-tens-products}
Let $M$ and $N$ be as in example \ref{ex_der-tens-products} and
set $P:=(\perp^R_\bullet M)\otimes_Rs.N$ (this is a simplicial
$R$-module; especially, it is a bisimplicial $A$-module); by
Eilenberg-Zilber's theorem \ref{th_Eilenberg-Zilber}, the cochain
complex associated with $M\ellotimes_RN$ is naturally isomorphic,
in $\sD(A\Mod)$ to the complex $\Tot(P^{\bullet\bullet})$, where
$P^{\bullet\bullet}$ denotes the double complex associated with $P$.
There follows a spectral sequence :
$$
E^1_{pq}:=\Tor_q^{R[p]}(M[p],N[p])\Rightarrow H_{p+q}(M\ellotimes_RN)
$$
whose differentials $d^1_{pq}:E^1_{pq}\to E^1_{p-1,q}$ are obtained as
follows. For every integer $q\in\N$, let $T_q$ be the $R$-module
given by the rule : $T_q[p]:=E^1_{pq}$, with face and degeneracy
maps deduced from those of $M$, $N$ and $R$, in the obvious way.
Then $d^1_{pq}$ is the differential in degree $p$ of the chain
complex $T_{q\bullet}$ associated with $T_q$ (details left to the
reader).
\end{remark}

\begin{example}\label{ex_standard-for-alg}
(i)\ \ 
Another useful triple arises from the forgetful functor
$A\Alg\to\Set$ and its left adjoint, that attaches to
any set $S$ the free $A$-algebra $A[S]$.
There results, for every $A$-algebra $B$, a standard
simplicial resolution by free $A$-algebras
$$
F_A(B)\to B.
$$
Now, if $R$ is again a simplicial $A$-algebra, and $S$ any
$R$-algebra, we can proceed as in the foregoing, to obtain
a bisimplicial resolution $F_R(S)\to S$ whose columns
$F_{R[n]}(S[n])\to S[n]$ are augmented simplicial $R[n]$-algebras,
for every $n\in\N$. A simple inspection shows that the proof
of lemma \ref{lem_Whitegead} carries over -- {\em mutatis
mutandis} -- to the present situation, hence a quasi-isomorphism
$S\to S'$ of $R$-algebra induces a morphism
$F_R(S)[n]\to F_R(S')[n]$ of $R$-algebras, which is a
homotopy equivalence on the underlying $R$-modules, for
every $n\in\N$.

(ii)\ \
We may then define a {\em derived tensor product} for $R$-algebras,
proceeding as in example \ref{ex_der-tens-products}. Namely,
if $S$ and $S'$ are any two $R$-algebras, we define two functors
$$
S\ellotimes_R-
\quad
\text{(\ resp. $-\ellotimes_RS'$\ )}
\quad : \quad
R\Alg\to R\Alg
$$
by the rules :
$$
T\mapsto S\otimes_RF_R(T)^\Delta
\quad
\text{(\ resp. $T\mapsto F_R(T)^\Delta\otimes_RN$\ )}
\qquad
\text{for every $R$-algebra $T$}.
$$
Arguing like in {\em loc.cit.} we see that both functors
transform quasi-isomorphisms into quasi-isomorphisms,
hence they descend to the derived category $\sD(R\Alg)$,
and moreover, the two functors are naturally isomorphic,
so the notation $S\ellotimes_RS'$ is unambiguous : the
detailed verification is left as an exercise for the reader.
\end{example}

\sset\subsubsection{}\label{subsec_new-subsection}
Let $A\tdu\AlgMod$ be the category of all pairs $(B,M)$,
where $B$ is any $A$-algebra, and $M$ any $B$-module.
The morphisms $(B,M)\to(B',M')$ are the pairs $(f,\phi)$,
where $f:B\to B'$ is a morphism of $A$-algebras, and
$\phi:M\to f^*M'$ a $B$-linear map (here $f^*M'$ denotes
the $B$-module obtained from the $B'$-module $M'$, by
restriction of scalars along the map $f$). Suppose now
that we have a commutative diagram of categories :
$$
\xymatrix{
A\tdu\AlgMod \ar[rr]^-T \ar[d]_F & & \cC \ar[d]^f \\
A\Alg \ar[rr]^-g & & \cB
}$$
such that :
\begin{itemize}
\item
For every $X\in\Ob(B)$, the fibre $f^{-1}X$ is an abelian
category whose filtered colimits are representable and exact.
\item
$F$ is given by the rule $(B,M)\mapsto B$ for every
$(B,M)\in\Ob(A\tdu\AlgMod)$.
\item
For every $A$-algebra $B$, the restriction
$$
T_B:F^{-1}B\to f^{-1}(gB)
$$
of $T$ commutes with filtered colimits, {\em i.e.}, if
$((B,M_i)~|~i\in I)$ is any filtered system of objects
of $A\tdu\AlgMod$ (over the same $A$-algebra $B$), then
the induced morphism
$$
\colim_{i\in I}T_BM_i\to T_B(\colim_{i\in I}M_i)
$$
is an isomorphism.
\end{itemize}
The functors $f$ and $g$ extend to functors $s.f:s.\cC\to s.\cB$,
respectively $s.g:s.A\Alg\to s.\cB$, and notice that -- for $R$
any simplicial $A$-algebra -- $T$ induces a functor
$$
T_R:R\Mod\to s.\cC_{/R}:=s.f^{-1}(s.g R)
\qquad
(M[n]~|~n\in\N)\mapsto(T_{R[n]}M[n]~|~n\in\N).
$$
For every $R$-module $M$, we set
$$
LT_R(M):=T_R(\perp^R_\bullet\!M)^\Delta.
$$
The augmentation $T_R(\perp^R_\bullet\!M)\to T_RM$ can be
regarded as a morphism of bisimplicial objects
$$
T_R(\perp^R_\bullet\!M)\to s.T_RM
$$
(the columns of $s.T_RM$ are constant $\cC$-simplicial
objects), whence a morphism in $s.\cC_{/R}$ :
\set\begin{equation}\label{eq_non-add-der-fctr}
LT_R(M)\to s.T_RM^\Delta=T_RM
\qquad
\text{for every $R$-module $M$}.
\end{equation}
In many applications, $s.\cC_{/R}$ will be also an abelian
category, but anyhow we can state the following :

\begin{proposition}\label{prop_derive-non-add}
With the notation of \eqref{subsec_new-subsection}, suppose
that $\cC$ is an abelian category. Then the following holds :
\begin{enumerate}
\item
If $M$ is a flat $R$-module, then \eqref{eq_non-add-der-fctr}
is a quasi-isomorphism.
\item
If\/ $\phi:M\to N$ is a quasi-isomorphism of $R$-modules, then
the induced map
$$
LT_R(M)\to LT_R(N)
$$
is a quasi-isomorphism.
\end{enumerate}
\end{proposition}
\begin{proof}(i): The assumption means that $M[n]$ is a flat
$R[n]$-module for every $n\in\N$. Denote by $C$ the unnormalized
double complex associated with $T_R(\perp^R_\bullet\!M)$, and by
$C^\Delta$ (resp. by $D$) the unnormalized complex associated
with $LT_R(M)$ (resp. to $T_RM$). By Eilenberg-Zilber's theorem
\ref{th_Eilenberg-Zilber}, we have a natural quasi-isomorphism
\set\begin{equation}\label{eq_alex-whitney}
C^\Delta\to\Tot(C)
\end{equation}
given by the Alexander-Whitney map of \eqref{subsec_AW-and-Sh},
and a simple inspection shows that the map $C^\Delta\to D$
induced by \eqref{eq_non-add-der-fctr} factors through
\eqref{eq_alex-whitney}. Hence, it suffices to show that
the map
$$
\Tot(C)\to D
$$
induced by the augmentation of $T_R(\perp^R_\bullet M)$, is
a quasi-isomorphism, under the stated condition.

To this aim, it suffices to check that the augmented
$\cC$-simplicial object
\set\begin{equation}\label{eq_short-name}
T_R(\perp^{R[n]}_\bullet\!M[n])\to T_RM[n]
\end{equation}
is aspherical for every $n\in\N$. However, since $M[n]$ is a
flat $R[n]$-module, it can be written as a filtered colimit of
a system of free $R[n]$-modules (\cite[Ch.I, Th.1.2]{La}); on
the other hand, the functor $\perp^{R[n]}_\bullet$ commutes with
all filtered colimits, and the same holds for $T_R$, by assumption.
Since the filtered colimits of the fibres of the functor $f$ are
exact, we are then reduced to checking the assertion in
case $M[n]$ is a free $R[n]$-module; but in this case,
\eqref{eq_short-name} is even homotopically trivial, by virtue
of proposition \ref{prop_triple-res}.

(ii): In light of Eilenberg-Zilber's theorem \ref{th_Eilenberg-Zilber},
it suffices to show that the induced map
$$
T_R(\perp_n^R\!M)\to T_R(\perp_n^R\!N)
$$
is a quasi-isomorphism for every $n\in\N$. In turns, this
follows readily from lemma \ref{lem_Whitegead}.
\end{proof}

\begin{remark}\label{rem_loop-and-suspend}
(i)\ \
The proof of proposition \ref{prop_derive-non-add}(i) applies
as well to the derived tensor product; namely, for any two
$R$-modules $M$ and $N$ there is a natural morphism of
$R$-modules
$$
M\ellotimes_RN\to M\otimes_RN
$$
which is a quasi-isomorphism, if either $M$ or $N$ is a
flat $R$-module (details left to the reader).

(ii)\ \
Especially, let $S$ and $S'$ be any two $R$-algebras.
Then, since the standard free resolution $F_R(S)$ of
example \ref{ex_standard-for-alg} is a flat simplicial
$R$-module, (i) implies that the $R$-module underlying
the $R$-algebra $S\ellotimes_RS'$, is naturally isomorphic,
in $\sD(R\Mod)$, to the derived tensor product of the
$R$-modules underlying $S$ and $S'$. In other words,
the notation $-\ellotimes_R-$ is unambiguous, whether
one refers to derived tensor products of algebras, or
of their underlying modules.

(iii)\ \
Denote by $\sigma,\omega:R\Mod\to R\Mod$ respectively
the suspension and loop functors (\cite[I.3.2.1]{Il}),
and recall that $\sigma$ is left adjoint to $\omega$.
A simple inspection of the definitions yields natural
identifications
$$
(\sigma M)\otimes_RN\isom\sigma(M\otimes_RN)
\qquad
(\omega M)\otimes_RN\isom\omega(M\otimes_RN)
$$
for every $R$-modules $M$ and $N$. By the same token,
it is clear that $\sigma$ and $\omega$ transform flat
$R$-modules into flat $R$-modules. In view of (i), we
deduce natural isomorphisms
$$
(\sigma M)\ellotimes_RN\isom\sigma(M\ellotimes_RN)
\qquad
(\omega M)\ellotimes_RN\isom\omega(M\ellotimes_RN)
\qquad
\text{in $\sD(R\Mod)$}
$$
for every $R$-modules $M$ and $N$.

(iv)\ \
Let $f:N\to N'$ be any morphism of $R$-modules; in the
same vein, we get a natural identification :
$$
\Cone(M\otimes_Rf)\isom M\otimes_R\Cone\,f
\qquad
\text{for every $R$-module $M$}
$$
(see \cite[I.3.2.2]{Il} for the definition of the cone of
a morphism of $R$-modules). It follows immediately that
the functor $M\ellotimes_R-:\sD(R\Mod)\to\sD(R\Mod)$ transforms
distinguished triangles into distinguished triangles (see
\cite[I.3.2.2.4]{Il} for the definition of distinguished
triangle in $\sD(R\Mod)$). For future reference, let us also
point out :
\end{remark}

\begin{lemma}\label{lem_usual-vanishing}
Let $R$ be a simplicial $A$-algebra, $X$ and $Y$ two
$R$-modules, $n,m\in\N$ two integers, and suppose that
$H_iX=0=H_jY$ for every $i<n$ and $j<m$. Then
$$
H_i(X\ellotimes_RY)=0
\qquad
\text{for every $i<n+m$.}
$$ 
\end{lemma}
\begin{proof} Set $\sigma^0:=\one_{R\Mod}$, and define
inductively $\sigma^k:R\Mod\to R\Mod$ by the rule :
$\sigma^k:=\sigma\circ\sigma^{k-1}$ for every $k>0$.
Define likewise $\omega^k$, and notice that $\sigma^k$
is left adjoint to $\omega^k$, for every $k\in\N$
(remark \ref{rem_loop-and-suspend}(iii)). Let
$f:\sigma^n\circ\omega^nX\to X$ denote the counit of
adjunction, and notice that $\Cone\,f=0$ in $\sD(R\Mod)$.
In view of remark \ref{rem_loop-and-suspend}(iv,v), we
deduce that the natural morphisms
$$
\sigma^n(\omega^nX\ellotimes_RY)\to
(\sigma^n\circ\omega^nX)\ellotimes_RY\to
X\ellotimes_RY
$$
are isomorphisms in $\sD(R\Mod)$. Likewise, the natural
morphism
$$
\sigma^m(\omega^nX\ellotimes_R\omega^mY)\to
\omega^nX\ellotimes_RY
$$
is an isomorphism in $\sD(R\Mod)$, so finally the
same holds for the morphism
$$
\sigma^{n+m}(\omega^nX\ellotimes_R\omega^mY)\to
X\ellotimes_RY
$$
whence the claim.
\end{proof}

\sset\subsubsection{}
In view of proposition \ref{prop_derive-non-add}, it is
clear that if $\cC$ is an abelian category, $LT_R$ yields
a well defined functor on derived categories
$$
LT_R:\sD(R\Mod)\to\sD(s.\cC_{/R})
\qquad
M\mapsto LT_R(M)
$$
(where $\sD(s.\cC_{/R})$ is the localization of $s.\cC_{/R}$
with respect to the class of quasi-isomorphisms, and likewise
for $\sD(R\Mod)$ : see \cite[Def.8.1.3]{Ga-Ra}). More
generally, suppose that $\cC$ is endowed with a functor
$$
\Phi:\cC\to\cA
$$
to an abelian category $\cA$ (usually, this will be a forgetful
functor of some sort), and denote by $\sD_\Phi(s.\cC_{/R})$ the
localization of $s.\cC_{/R}$ with respect to the system of its
morphisms whose image under $s.\Phi$ are quasi-isomorphisms in
$s.\cA$. Then, since clearly
$$
s.\Phi\circ LT_R=L(s.\Phi\circ T_R)
$$
we see that $LT_R$ descends to a well defined functor
$$
LT_R:\sD(R\Mod)\to\sD_\Phi(s.\cC_{/R}).
$$
We will need also the following slight refinement :

\begin{corollary}\label{cor_truncate-and-derive}
In the situation of proposition {\em\ref{prop_derive-non-add}},
the following holds :
\begin{enumerate}
\item
Let $\phi:M\to N$ be a morphism of $R$-modules, and $n\in\N$
an integer such that
$$
H_i\phi:H_iM\to H_iN
$$
is an isomorphism, for every $i<n$. Then the same holds for
the induced map
$$
H_i(LT_R\phi):H_i(LT_RM)\to H_i(LT_RN).
$$
\item
Suppose that $T(B,0)=0_{gB}$ (the initial and final object of
$f^{-1}(gB)$) for every $A$-algebra $B$. Let $n\in\N$ be an integer,
and $M$ an $R$-module such that $H_iM=0$ for every $i<n$.
Then $H_i(LT_RM)=0_{gR_i}$ for every $i<n$.
\end{enumerate}
\end{corollary}
\begin{proof}(i): Denote by
$$
\phi_X:X\to\cosk_nX
\qquad
\text{for every $X\in\Ob(R\Mod)$}
$$
the unit of adjunction (see \eqref{subsec_coskeleton}). Taking
into account corollary \ref{cor_Dold-Kan}(i), we have the
following properties :
\begin{itemize}
\item
$X[i]=\cosk_nX[i]$ and $\phi[i]$ is the identity map of
$X[i]$, for every $i\leq n$.
\item
$H_i(\cosk_nX)=0$ for every $i\geq n$.
\end{itemize}
In view of proposition \ref{prop_derive-non-add}, we are then
easily reduced to checking the assertion for the special where
$N:=\cosk_nM$, and $\phi:=\phi_M$.
However, since in this case $N[i]=M[i]$ for every $i\leq n$,
it is clear that $\perp^{R[i]}_\bullet\!\phi[i]$ is the identity
automorphism of $\perp_\bullet^{R[i]}\!\!M[i]$, for every $i\leq n$.
After applying the functor $T_R$, and extracting the diagonal,
we see that the induced morphism $LT_RM\to LT_RN$ in $s.\cC$ is
given by the identity automorphism of $(LT_RM)[i]$, in every
degree $i\leq n$, whence the contention.

(ii): In light of (i), we may assume that $M=s.0$ (the trivial
$R$-module). In this case, the assertion follows easily from
remark \ref{rem_trivial-to-trival}.
\end{proof}

\sset\subsubsection{}\label{subsec_T-and-base-change}
In the situation of \eqref{subsec_new-subsection}, take
$\cC:=A\tdu\AlgMod$, $f:=F$, and $g$ the identity
automorphism of $A\Alg$. If $B\to B'$ is a map of
$A$-algebras, and $(B,M)\in\Ob(\cC)$, we define
$B'\otimes_B(B,M):=(B',B'\otimes_BM)$. 

\begin{corollary}
In the situation of \eqref{subsec_T-and-base-change},
suppose moreover that, for every flat morphism $B\to B'$
of\/ $A$-algebras, the natural map
$$
B'\otimes_BT(B,M)\to T(B'\otimes_B(B,M))
$$
is an isomorphism, for every $(B,M)\in\Ob(\cC)$. Then,
for every flat morphism $\phi:R\to S$ of simplicial
$V$-algebras, the natural morphism
\set\begin{equation}\label{eq_natural-map}
S\otimes_RLT_RM\to LT_S(S\otimes_RM)
\end{equation}
is an isomorphism in $\sD(S\Mod)$, for every $R$-module $M$.
\end{corollary}
\begin{proof} (Recall that $\phi$ is flat if and only if
$S[n]$ is a flat $R[n]$-algebra, for every $n\in\N$. The
map of the proposition is deduced from the natural morphism
$\perp^R_\bullet\!M\to\perp^S_\bullet\!(S\otimes_RM)$, given by
functoriality of the standard free resolution.)

Let $\pi:\perp^R_\bullet\!M\to M$ be the standard free resolution
of $M$; by the flatness assumption on $S$, the morphism
$S\otimes_R\pi:S\otimes_R\perp^R_\bullet\!M\to S\otimes_RM$
is a free resolution of the $S$-module $S\otimes_RM$, whence
a natural isomorphism
$$
LT_S(S\otimes_RM)\isom T_S((S\otimes_R\perp^R_\bullet\!M)^\Delta)
\qquad
\text{in $\sD(S\Mod)$}
$$
by proposition \ref{prop_derive-non-add}(i) and
Eilenberg-Zilber's theorem \ref{th_Eilenberg-Zilber}. On the
other hand, the assumptions on $T$ yield a natural isomorphism
$$
S\otimes_RLT_R(M)=
S\otimes_RT_R((\perp^R_\bullet\!M)^\Delta)\isom
T_S((S\otimes_R\perp^R_\bullet\!M)^\Delta)
\qquad
\text{in $\sD(S\Mod)$}.
$$
By combining these isomorphisms, we get an isomorphism
$S\otimes_RLT_RM\to LT_S(S\otimes_RM)$, and a simple inspection
shows that the latter is realized by the natural map
\eqref{eq_natural-map}.
\end{proof}

\sset\subsubsection{}\label{subsec_quasi-regular}
Let now $R$ be any simplicial $A$-algebra, and $I\subset R$
any ideal, and set $R_0:=R/I$. The {\em Rees algebra of $I$}
is the graded $R$-algebra
$$
\sR(R,I)^\bullet:=\bigoplus_{p\in\N}I^p
$$
(cp. example \ref{ex_Rees}). Notice the natural isomorphism
of graded $R_0$-algebras
$$
\sR(R,I)\otimes_RR_0\isom\gr^\bullet_IR:=
\bigoplus_{p\in\N}I^p/I^{p+1}
$$
as well as the exact sequence of $\sR(R,I)$-modules
$$
0\to\gr^{\bullet+1}_IR\to\sR(R,I)\otimes_RR/I^2\to
\gr^\bullet_IR\to 0
$$
deduced from the natural projection $R/I^2\to R_0$. Next,
pick any quasi-isomorphism of $R$-algebras $P\to R_0$,
with $P$ a flat $R$-algebra; there follows an exact sequence
of $P\otimes_R\sR(R,I)$-modules
$$
0\to P\otimes_R\gr^{\bullet+1}_IR\to P\otimes_R\sR(R,I)\otimes_RR/I^2
\xrightarrow{\ \beta\ }P\otimes_R\gr^\bullet_IR\to 0
$$
and notice that $\beta$ is actually a morphism of
$P\otimes_R\sR(R,I)$-algebras. According to remark
\ref{rem_from-simpl-to-dga}(i), after forming the associated
cochain complexes, we obtain an epimorphism $\beta^\bullet$
of differential graded $(P\otimes_R\sR(R,I))^\bullet$-algebras,
whose kernel $(P\otimes_R\gr^{\bullet+1}_IR)^\bullet$ is a
two-sided ideal of $(P\otimes_R\sR(R,I)\otimes_RR/I^2)^\bullet$.
In this situation, we deduce a map of
$H_\bullet(R_0\ellotimes_R\sR(P,I))$-modules
\set\begin{equation}\label{eq_here-a-deriv}
\delta:H_\bullet(R_0\ellotimes_R\gr^\bullet_IR)\to
H_\bullet(R_0\ellotimes_R\gr^{\bullet+1}_IR)[1]
\end{equation}
which is a graded derivation, according to lemma
\ref{lem_finally-got-it} (recall that here the shift $[1]$
refers to the homological grading; the notation $\gr^{\bullet+1}$
denotes also a shift in degrees, which however {\em does not}
alter the signs of the scalar multiplication in the way
prescribed by remark \ref{rem_dgas}(v)).

\sset\subsubsection{}\label{subsec_quasi-regular-true}
Keep the situation of \eqref{subsec_quasi-regular}, and
suppose furthermore that $R$ is a constant $A$-algebra
and $I$ a constant ideal, {\em i.e.} $R=s.B$, and $I=s.J$
for some $A$-algebra $B$ and some ideal $J\subset B$.
We set $B_0:=B/J$ and
$$
G^\bullet:=R_0\ellotimes_R\gr_I^\bullet R
\qquad
\Lambda_\bullet:=\Lambda^\bullet_B(J/J^2)
\qquad
S^\bullet:=\Sym_B^\bullet(J/J^2)
$$
where the derived tensor product is formed in $\sD(R\Alg)$,
as in example \ref{ex_standard-for-alg}, so $G^\bullet$ is
represented by $P\otimes_R\gr^\bullet_IR$, and $H_\bullet G^\bullet$
is a bigraded $A$-algebra, strictly anti-commutative for the
(homological) grading $H_\bullet$, and commutative for
the (cohomological) grading deduced from the grading of
$G^\bullet$ (see remark \ref{rem_from-simpl-to-dga}(i)).
Also, $\Lambda_\bullet$ (resp. $S^\bullet$) is a strictly
anti-commutative (resp. a commutative) graded $B_0$-algebra,
and we use the homological (resp. cohomological) conventions
for grading, so $\Lambda_p$ (resp. $S^p$) is placed in
degree $-p$ (resp. $p$).
Moreover, a standard calculation yields a natural isomorphism
of $B_0$-modules
$$
H_0G^0\isom B_0
\qquad
H_1G^0\isom\Tor_1^B(B_0,B_0)\isom J/J^2
$$
whence a natural map $\Lambda_\bullet\to H_\bullet G^0$ of
strictly anti-commutative graded $B_0$-algebras, restricting
to an isomorphism in degrees $\leq 1$. On the other hand, we
have a surjective map of commutative graded $B_0$-algebras
\set\begin{equation}\label{eq_symm-powers}
S^\bullet\to\gr^\bullet_JB:=\bigoplus_{p\in\N}J^p/J^{p+1}
\end{equation}
restricting as well to an isomorphism in degrees $\leq 1$;
since $P\otimes_R\gr^\bullet_IR$ is naturally a simplicial
$\gr^\bullet_JB$-algebra, we obtain a natural map of bigraded
$S^\bullet$-algebras
$$
\omega^\bullet_\bullet:
\Lambda_\bullet\otimes_BS^\bullet\to H_\bullet G^\bullet.
$$
Now, $H_\bullet G^\bullet$ has been endowed with a derivation
$\delta:H_\bullet G^\bullet\to H_\bullet G^{\bullet+1}[1]$
in \eqref{eq_here-a-deriv}; on the other hand, there exists
on $\Lambda_\bullet\otimes_BS^\bullet$ a natural $S^\bullet$-linear
graded derivation
$$
\partial:
\Lambda_\bullet\otimes_BS^\bullet\to
\Lambda_\bullet\otimes_BS^{\bullet+1}[1]
$$
as explained in \cite[I.4.3.1.2]{Il}. We may then consider
the diagram of $S^\bullet$-linear maps
\set\begin{equation}\label{eq_partials}
{\diagram
\Lambda_\bullet\otimes_BS^\bullet \ar[rr]^-{\omega^\bullet_\bullet}
\ar[d]_\partial & & H_\bullet G^\bullet \ar[d]^\delta \\
\Lambda_\bullet\otimes_BS^{\bullet+1}[1]
\ar[rr]^-{\omega^{\bullet+1}_\bullet[1]} & & H_\bullet G^{\bullet+1}[1].
\enddiagram}
\end{equation}

\begin{definition}\label{def_Quillen}
(i)\ \ 
In the situation of \eqref{subsec_quasi-regular-true}, we say
that the ideal $J$ is {\em quasi-regular}, if the following
conditions hold :
\begin{itemize}
\item
The map $\omega^\bullet_\bullet$ restricts to an isomorphism
$\omega^0_\bullet\isom H_\bullet G^0$.
\item
The $B_0$-module $J/J^2$ is flat.
\end{itemize}

(ii)\ \
Let $R$ be any simplicial $A$-algebra, and $I\subset R$
any ideal. We say that $I$ is {\em quasi-regular}, if $I[n]$
is a quasi-regular ideal of $R[n]$, for every $n\in\N$.
\end{definition}

\begin{remark}\label{rem_augm-of-free}
(i)\ \
Definition \ref{def_Quillen} is due to Quillen (\cite[Def.8.4]{Qu0}).
Notice that a sequence of elements $(f_1,\dots,f_r)$ in a ring $B$
is $B$-quasi-regular, in the sense of definition
\ref{def_quasi-regular-seq}, if and only if it generates an ideal
$J\subset B$ such that $J/J^2$ is a free $B_0$-module of rank $r$,
and \eqref{eq_symm-powers} is an isomorphism (notation of
\eqref{subsec_quasi-regular-true}). The relationship between these
two quasi-regularity notions is partially elucidated by the
following proposition \ref{prop_quasi-regular}(ii).

(ii)\ \ 
Suppose that $\underline B:=(B_i~|~i\in I)$ is any filtered
system of rings, and $\underline J:=(J_i\subset B_i~|~i\in I)$
a filtered system of ideals, such that $J_i$ is quasi-regular
in $B_i$, for every $i\in I$. Denote by $B$ (resp. $J$) the
colimit of $\underline B$ (resp. of $\underline J$); then
it is easily seen that $J$ is quasi-regular in $B$.

(iii)\ \
Let $B$ be any ring, and $C:=B[X_i~|~i\in I]$ a free $B$-algebra
(for any set $I$). Then the ideal $J:=(X_i~|~i\in I)\subset C$ is
quasi-regular. Indeed, (ii) reduces to the case where $I$ is a
finite set, for which the assertion is a special case of
proposition \ref{prop_quasi-reg-equals-reg} below.
\end{remark}

\begin{proposition}\label{prop_quasi-regular}
In the situation of \eqref{subsec_quasi-regular-true}, we have :
\begin{enumerate}
\item
\eqref{eq_partials} is a commutative diagram.
\item
Suppose that $J$ is a quasi-regular ideal. Then
\begin{enumerate}
\item
$\omega^\bullet_\bullet$ is an isomorphism.
\item
There exists a natural isomorphism of $B_0$-modules :
$$
\Tor_i^B(B_0,B/J^n)\isom\Coker
(\partial^{n-2}_{i+1}:\Lambda_{i+1}\otimes_BS^{n-2}\to
\Lambda_i\otimes_BS^{n-1})
\qquad
\text{for every $n,i>0$}.
$$
\end{enumerate}
\end{enumerate}
\end{proposition}
\begin{proof}(i): Since both derivations are $S^\bullet$-linear,
it suffices to check that \eqref{eq_partials} commutes in
(upper) degree $0$. Morever, since the $B$-algebra
$\Lambda_\bullet$ is generated by $\Lambda_1$, it suffices
to check the commutativity of the diagram
\set\begin{equation}\label{eq_partials_10}
{\diagram
\Lambda_1 \ar[r]^{\omega^0_1} \ar[d]_-{\partial_1^0} &
H_1G^0 \ar[d]^{\delta^0_1} \\
S^1 \ar[r]^-{\omega_0^1} & H_0G^1.
\enddiagram}
\end{equation}
However, according to \cite[I.4.3.1.2]{Il}, the map
$\partial^0_1$ is just the identity of $J/J^2$.
The map $\delta^0_1$ is the boundary map arising -- by the
snake lemma -- from the exact sequence of complexes
$$
0\to JP^\bullet/J^2P^\bullet\to P^\bullet/J^2P^\bullet\to
P^\bullet/JP^\bullet\to 0
$$
(where $P^\bullet$ denotes the cochain complex associated
with the flat resolution $P\to s.B_0$). The same boundary
map is used to define the isomorphism $\omega_1^0$, and
$\omega^1_0$ is the natural isomorphism
$J/J^2\to H_0(P^\bullet/JP^\bullet)$. Summing up, the
commutativity of \eqref{eq_partials_10} follows by simple
inspection.

(ii.a): We prove, by induction on $n\in\N$, that every
$\omega_\bullet^n$ is an isomorphism. For $n=0,1$ there
is nothing to show, so suppose that $n>1$, and the
assertion is already known for every degree $<n$.
In order to prove the assertion in degree $n$, it suffices
to check that $\omega^n_0:S^n\to\gr^n_JB$ is an isomorphism;
indeed, in this case $\gr^n_JB$ is a flat $B_0$-module,
and then a standard spectral sequence argument allows to
conclude that $\omega^n_j$ is also an isomorphism, for every
$j\in\N$. 

Let $P\to R_0$ be a resolution as in \eqref{subsec_quasi-regular}.
The $I$-adic filtration has finite length on $P/I^kP$, for every
$k\in\N$, and therefore yields a convergent spectral sequence
$$
E^1_{pq}\Rightarrow \Tor_{p+q}^B(B_0,B/J^k)
$$
with
$$
E^1_{pq}:=\left\{\begin{array}{ll}
           \Tor^B_{p+q}(B_0,\gr_JB^{-p}) & \text{for $p>-k$} \\
           0 & \text{otherwise}
          \end{array}\right.
$$
whose differentials $E^1_{pq}\to E^1_{p-1,q}$ agree -- by direct
inspection -- with the derivation $\delta^{-p}_{p+q}$, whenever
$p>1-k$. Especially, for $k=n+1$, we get a complex
$$
\Sigma
\quad :\quad
H_2G^{n-2}\xrightarrow{\ \delta_2^{n-2}\ }
H_1G^{n-1}\xrightarrow{\ \delta_1^{n-1}\ }
\gr^n_JB\to 0.
$$
On the other hand, since $J/J^2$ is a flat $B_0$-module,
the corresponding sequence 
$$
\Sigma'
\quad :\quad
\Lambda_2\otimes_BS^{n-2}\xrightarrow{\ \partial\ }
\Lambda_1\otimes_BS^{n-1}\xrightarrow{\ \partial\ }S^n\to 0
$$
is an exact complex (this is a segment of the Koszul complex
of \cite[I.4.3.1.7]{Il}). In light of (i), the maps
$\omega^{n-2}_2$ and $\omega^{n-1}_1$ yield a morphism of
complexes $\Sigma'\to\Sigma$, so we are reduced to checking
that $\Sigma$ is exact, and the inductive assumption already
says that $\Sigma$ is exact at the middle term $H_1G^{n-1}$.
Now notice that, since the Koszul complex is exact, the
inductive assumption implies that $E^2_{pq}=0$ for all $(p,q)$
such that $p>1-n$ and $q>0$; moreover, obviously we have
$E^1_{pq}=0$ for $p+q<0$ or $p>0$, therefore
$$
E^\infty_{00}=E^1_{00}
\qquad\text{and}\qquad
E^\infty_{-n,n}=E^2_{-n,n}=\Coker\,\delta^{n-1}_1.
$$
But a simple inspection shows that the natural map
$H_0P\to E^\infty_{00}$ is an isomorphism, so
$E^\infty_{-n,n}=0$, therefore $\delta^{n-1}_1$ is
surjective, and $\Sigma$ is exact, as required.

(ii.b): We use the previous spectral sequence, for $k=n$.
Taking into account (ii.a), and the exactness of the
Koszul complex (\cite[I.4.3.1.7]{Il}) we get
$$
\begin{aligned}
& E^1_{pq}=0 & \qquad & \text{for $p+q<0$ or $p>0$ or $p\leq -n$} \\
& E^2_{p,q}=0 & \qquad &
\text{for every $(p,q)$ with $p>1-n$ and $q>0$} \\
& E^2_{1-n,i+n-1}\isom\Coker\,\partial^{n-2}_{i+1} & \qquad &
\text{for every $i\in\N$}.
\end{aligned}
$$
It follows that $E^r_{1-n+r,i+n-r}=0$ for every $r\geq 2$ and
every $i\in\N$; indeed, this is clear if $1-n+r>0$, and if
the latter fails, we get $1-n+r>1-n$ and $i+n-r>0$, so the
stated vanishing holds nevertheless. We deduce that
$E^\infty_{1-n,i+n-1}=E^2_{1-n,i+n-1}$ for every $i\in\N$.
Next, suppose that $i>0$. Clearly we have $E^\infty_{p,i-p}=0$
both for $p>0$ and $p\leq -n$; and if $0\geq p>1-n$, we get
$i-p\geq i>0$, so again $E^\infty_{p,i-p}=0$. In conclusion,
$E^\infty_{p,i-p}=0$ except possibly for $p=1-n$, and the
contention follows.
\end{proof}

\begin{proposition}\label{prop_quasi-reg-equals-reg}
Let $B$ be a ring, $f_\bullet:=(f_1,\dots,f_r)$
a finite sequence of elements of $B$. Denote by $J$
the ideal generated by $f_\bullet$, set $B_0:=B/J$,
and consider the following conditions :
\begin{enumerate}
\alphaenu
\item
The sequence $f_\bullet$ is $B$-regular.
\item
The sequence $f_\bullet$ is completely secant in $B$.
\item
The ideal $J$ is quasi-regular and $J/J^2$ is a
free $B_0$-module of rank $r$.
\item
The sequence $f_\bullet$ is $B$-quasi-regular.
\end{enumerate}
Then {\em (a)$\Rightarrow$(b)$\Rightarrow$(c)$\Rightarrow$(d)},
and if $B$ is $J$-adically complete and separated,
then {\em (d)$\Rightarrow$(a)}.
\end{proposition}
\begin{proof}(a)$\Rightarrow$(b) and (c)$\Rightarrow$(d) have
already been pointed out respectively in proposition
\ref{prop_Kosz-cptl-sec} and proposition
\ref{prop_quasi-regular}(ii.a), and the same for the implication
(d)$\Rightarrow$(a), when $B$ is $J$-adically complete and
separated.

(b)$\Rightarrow$(c): This becomes a paraphrase of proposition
\ref{prop_Koszul-is-back}, once we have identified the graded
$B$-algebra $H_\bullet G^\bullet$ with the graded $B$-algebra
$H_\bullet(B_0[0]\derotimes_BB_0[0])$. However, denote by
$P^\bullet$ the cochain complex associated with the flat
resolution $P\to s.B_0$; according to remark
\ref{rem_from-simpl-to-dga}(i), $P^\bullet$ is naturally
a differential graded $B$-algebra, and the induced map
$\eps^\bullet:P^\bullet\to B_0[0]$ is a morphism of differential
graded $B$-algebras (for the multiplication law
$\bar\mu{}^\bullet$ on $B_0[0]$ inherited from that of $B_0$).
Then there exists in $\sD(B\Mod)$ a unique isomorphism
$\omega^\bullet:P^\bullet_{B_0}\to P^\bullet$ whose composition
with $\eps^\bullet$ agrees with $\rho^\bullet_{B_0}$ (notation
of \eqref{subsec_give_Ps}), and moreover we get a commutative
diagram in $\sD(B\Mod)$ :
$$
\xymatrix{ P^\bullet_{B_0}\otimes_BP^\bullet_{B_0}
\ar[rr]^-{P^\bullet_{\bar\mu}}
\ar[d]_{\omega^\bullet\otimes_B\omega^\bullet} & &
P^\bullet_{B_0} \ar[d]^{\omega^\bullet} \\
P^\bullet\otimes_BP^\bullet \ar[rr] & &
P^\bullet
}$$
whose bottom horizontal arrow is the multiplication map
of $P^\bullet$, and where $P^\bullet_{\bar\mu}$ is defined
as in \eqref{subsec_give_Ps}. The assertion follows
immediately.
\end{proof}

The following is a well known theorem of Quillen, which
is found in his typewritten notes \cite[Th.8.8]{Qu0} that
have been circulated since the late 60s, but have never
been published.

\begin{theorem}\label{th_Quillen}
Let $R$ be any simplicial $A$-algebra, and $I\subset R$ a
quasi-regular ideal such that $H_0I=0$. Then we have :
$$
H_iI^n=0
\qquad
\text{for every $n\in\N$ and every $i<n$}.
$$
\end{theorem}
\begin{proof} We argue by induction on $n\in\N$. For
$n\leq 1$ there is nothing to prove, so suppose that $n>1$,
and that the assertion is already known for $I^{n-1}$.
Set $R_k:=R/I^{k+1}$ for every $k\in\N$, and let
$$
\Lambda_\bullet\otimes_RS^\bullet:=
\Lambda^\bullet_R(I/I^2)\otimes_R\Sym^\bullet_R(I/I^2)
$$
which is a bigraded $R_0$-module, endowed with the
bigraded derivation $\partial^\bullet_\bullet$ such
that $\partial^\bullet_\bullet[i]$ is the Koszul
derivation of $\Lambda_\bullet\otimes_RS^\bullet[i]$,
defined as in \cite[I.4.3.1.2]{Il}. Set as well
$$
C_i:=\Coker\,(\partial^{n-3}_i:
\Lambda_i\otimes_RS^{n-3}\to\Lambda_{i-1}\otimes_RS^{n-2})
\qquad
\text{for every $i\in\N$}.
$$
From the inductive assumption and lemma
\ref{lem_usual-vanishing}, we already know that
\set\begin{equation}\label{eq_already-known-vanish}
H_i(I\ellotimes_RI^{n-1})=0
\qquad
\text{for every $i<n$}.
\end{equation}
However, according to remark \ref{rem_der-tens-products},
there is a spectral sequence
$$
E^1_{pq}:=\Tor^{R[p]}_q(I[p],I[p]^{n-1})\Rightarrow
H_{p+q}(I\ellotimes_RI^{n-1})
$$
and proposition \ref{prop_quasi-regular}(ii.b) yields natural
isomorphisms
$$
E^1_{pq}\isom\Tor^{R[p]}_{q+1}(R_0[p],I[p]^{n-1})
\isom\Tor^{R[p]}_{q+2}(R_0[p],R_{n-2}[p])
\isom C_{q+3}[p]
\qquad
\text{for every $q>0$}
$$
under which, the differentials $E^1_{pq}\to E^1_{p-1,q}$ are
identified with those of the chain complex associated with the
$R_0$-module $C_{q+3}$. Now, since $I/I^2$ is a flat
$R_0$-module, \cite[I.4.3.1.7]{Il} says that the sequence
of $R_0$-modules
$$
\Sigma^i
\quad : \quad
0\to\Lambda_{n+i-3}\xrightarrow{\ \partial\ }
\Lambda_{n+i-4}\otimes_RS^1\to\cdots\to
\Lambda_i\otimes_RS^{n-3}\xrightarrow{\ \partial\ }
\Lambda_{i-1}\otimes_RS^{n-2}\to C_i\to 0
$$
is an exact complex, for every $i>0$. On the other hand, since
$H_0I/I^2=0$, (and again, since $I/I^2$ is a flat $R_0$-module),
combining \cite[I.4.3.2.1]{Il} and lemma \ref{lem_usual-vanishing},
we see that
$$
H_i(\Lambda_j\otimes_RS^{q-j})=0
\qquad
\text{for every $q\in\N$, every $j\leq q$, and every $i<q$}.
$$
We recall the following general

\begin{claim}\label{cl_systematic}
Let $\cA$ be any abelian category, $d\in\N$ any integer, and
consider an exact sequence
$$
0\to K_d\xrightarrow{\ f_d\ } K_{d-1}\to\cdots\to
K_1\xrightarrow{\ f_1\ } K_0\to 0
$$
of objects of $\sC(\cA)$, such that $K_i=0$ in $\sD(\cA)$,
for every $i=1,\dots,d-1$. Then there is a natural isomorphism
$$
K_0\isom K_d[d-1]
\qquad
\text{in $\sD(\cA)$}.
$$
\end{claim}
\begin{pfclaim} For every $i=0,\dots,d$, set $Z_i:=\Ker\,f_i$.
We have short exact sequences
$$
0\to Z_i\to K_i\to Z_{i-1}\to 0
\qquad
\text{for $i=1,\dots,d$}.
$$
Since $K_i=0$ in $\sD(\cA)$, there follows natural isomorphisms
$Z_{i-1}\isom Z_i[1]$ for every $i=1,\dots,d-1$. However, $Z_0=K_0$
and $Z_{d-1}=K_d$, whence the claim.
\end{pfclaim}

By applying claim \ref{cl_systematic} to the exact sequence
$\cosk_{n+i-3}\Sigma^i$, we easily deduce that
$$
H_jC_i=0
\qquad
\text{for every $i>0$ and every $j<n+i-3$}
$$
so finally $E^2_{pq}=0$ for every $q>0$ and every $p<n+q$,
and clearly $E^1_{pq}=0$ for every $q<0$. Consequently,
for every $i<n$ the term $E^\infty_{i-q,q}=0$ vanishes
for $q>0$, and equals $E^2_{i,0}$ for $q=0$.
Lastly, a simple calculation shows that
$$
E^2_{i,0}=H_i(I\otimes_RI^{n-1})
\qquad
\text{for every $i\in\N$}.
$$
Summing up, and taking into account \ref{eq_already-known-vanish},
we conclude that $H_i(I\otimes_RI^{n-1})=0$ for every $i<n$,
so we are reduced to checking :

\begin{claim} The natural map
$H_i(I\otimes_RI^{n-1})\to H_i(I^n)$ is surjective,
for every $i\leq n$.
\end{claim}
\begin{pfclaim}[] Indeed, denote by $K$ the kernel of the
epimorphism $I\otimes_RI^{n-1}\to I^n$; as in the foregoing,
we get natural isomorphisms
$$
K[i]\isom\Tor_1^{R[i]}(R_0[i],I[i]^{n-1})\isom
\Tor_2^{R[i]}(R_0[i],R[i]/I[i]^{n-1})\isom
C_3[i]
\qquad
\text{for every $i\in\N$}
$$
that amount to an isomorphism $K\isom C_3$ of $R_0$-modules.
We have already seen that $H_iC_3=0$ for every $i<n$, whence
the claim.
\end{pfclaim}
\end{proof}

\subsection{Injective modules, flat modules and indecomposable
modules}

\sset\subsubsection{Indecomposable modules}
Recall that a unital (not necessarily commutative) ring $R$
is said to be {\em local\/} if $R\neq 0$ and, for every
$x\in R$ either $x$ or $1-x$ is invertible. If $R$ is
commutative, this definition is equivalent to the usual one.

\sset\subsubsection{}\label{subsec_indecomp}
Let $\cC$ be any abelian category and $M$ an object of $\cC$.
One says that $M$ is {\em indecomposable\/} if it is non-zero
and cannot be presented in the form $M=N_1\oplus N_2$ with
non-zero objects $N_1$ and $N_2$. If such a decomposition exists,
then $\End_\cC(N_1)\times\End_\cC(N_2)\subset\End_\cC(M)$, especially
the unital ring $\End_\cC(M)$ contains an idempotent element
$e\neq 1,0$ and therefore it is not a local ring.

However, if $M$ is indecomposable, it does not necessarily
follow that $\End_\cC(M)$ is a local ring. Nevertheless, one
has the following:

\begin{theorem}[Krull-Remak-Schmidt]\label{th_KRS}
Let $(A_i)_{i\in I}$ and $(B_j)_{j\in J}$ be two finite families
of objects of\/ $\cC$, such that:
\begin{enumerate}
\alphaenu
\item
$\oplus_{i\in I}A_i\simeq\oplus_{j\in J}B_j$.
\item
$\End_\cC(A_i)$ is a local ring for every $i\in I$, and
$B_j\neq 0$ for every $j\in J$.
\romanenu
\end{enumerate}
Then we have :
\begin{enumerate}
\item
There is a surjection $\phi:I\to J$, and isomorphisms
$B_j\isom\oplus_{i\in\phi^{-1}(j)} A_i$, for every $j\in J$.
\item
Especially, if $B_j$ is indecomposable for every $j\in J$,
then $I$ and $J$ have the same cardinality, and $\phi$ is
a bijection.
\end{enumerate}
\end{theorem}
\begin{proof} Clearly, (i)$\Rightarrow$(ii).
To show (i), let us begin with the following :

\begin{claim}\label{cl_immediate}
Let $M_1,M_2$ be two objects of $\cC$, and set $M:=M_1\oplus M_2$.
Denote by $e_i:M_i\to M$ (resp. $p_i:M\to M_i$) the natural injection
(resp. projection) for $i=1,2$. Suppose that $\alpha:P\to M$ is a
subobject of $M$, such that $p_1\alpha:P\to M_1$ is an isomorphism.
Then the natural morphism $\beta:P\oplus M_2\to M$ is an isomorphism.
\end{claim}
\begin{pfclaim} Denote by $\pi_P:P\oplus M_2\to P$ and
$\pi_2:P\oplus M_2\to M_2$ the natural projections; then
$\beta:=\alpha\pi_P+e_2\pi_2$. Of course, the assertion
follows easily by applying the $5$-lemma (which holds in any
abelian category) to the commutative ladder with exact rows :
$$
\xymatrix{
0\ar[r] & M_2 \ar[r] \ddouble &
P\oplus M_2 \ar[r]^-{\pi_P} \ar[d]_\beta &
P \ar[d]^{p_1\alpha} \ar[r] & 0 \\
0 \ar[r] & M_2 \ar[r]^-{e_2} & M \ar[r]^-{p_1} & M_1 \ar[r] & 0.
}$$
Equivalently, one can argue directly as follows. By definition,
$\Coker\,\beta$ represents the functor
$$
\cC\to\Z\Mod
\quad :\quad
X\mapsto\Ker\,\Hom_\cC(\beta,X)
$$
and likewise for $\Coker\,p_1\alpha$; however, it is easily
seen that these functors are naturally isomorphic, hence the
natural morphism $\Coker\,\beta\to\Coker\,p_1\alpha$ is an
isomorphism, so $\beta$ is an epimorphism. Next, let
$t:\Ker\,\beta\to P\oplus M_2$ be the natural morphism;
since $p_1\alpha$ is an isomorphism, $\pi_P\circ t=0$,
so $t$ factors through a morphism $\bar t:\Ker\,\beta\to M_2$.
It follows that $e_2\circ\bar t=\beta\circ t=0$, therefore
$\bar t=0$, and finally $\Ker\,\beta=0$, since $t$ is a
monomorphism.
\end{pfclaim}

Set $A:=\oplus_{i\in I} A_i$, and denote $e_i:A_i\to A$ (resp.
$p_i:A\to A_i$) the natural injection (resp. projection), for
every $i\in I$. Also, endow $I':=I\cup\{\infty\}$ with any total
ordering $(I',<)$ for which $\infty$ is the maximal element.

\begin{claim}\label{cl_a-plus-b}
Let $(a_j:A\to A)_{j\in J}$ be a finite system of endomorphisms
such that $\sum_{j\in J}a_j=\one_A$. Then there exists a mapping
$\phi:I\to J$ such that the following holds :
\begin{enumerate}
\item
$a_{\phi(i)}e_i:A_i\to A$ is a monomorphism, for every $i\in I$,
and let $P_i:=\Img(a_{\phi(i)}e_i)$.
\item
For every $l\in I'$, the natural morphism
$\beta_l:\bigoplus_{i<l} P_i\oplus\bigoplus_{i\geq l}A_i\to A$
is an isomorphism.
\end{enumerate}
\end{claim}
\begin{pfclaim} We both define $\phi(l)$ and prove (ii), by
induction on $l$. If $l$ is the smallest element of $I$, then
$\beta_l=\one_A$, so (ii) is obvious. Hence, let $l\in I$ be any
element, $l'\in I'$ the successor of $l$, and assume that $\beta_l$
is an isomorphism. Set $a'_j:=\beta_l^{-1}a_j\beta_l$ for every
$j\in J$, and $M_l:=\bigoplus_{i<l} P_i\oplus\bigoplus_{i\geq l}A_i$;
also, let $e'_l:A_l\to M_l$ and $p'_l:M_l\to A_l$ be the natural
morphisms. Clearly $\sum_{j\in J}a'_j=\one_{M_l}$, so
$\sum_{j\in J}p'_la'_je'_l=\one_{A_l}$; since $\End_\cC(A_l)$
is a local ring, it follows that there exists $j_l\in J$ such
that $p'_la'_{j_l}e'_l$ is an automorphism of $A_l$. However,
\set\begin{equation}\label{eq_Krull-Remark}
a'_{j_l}e'_l=\beta_l^{-1}a_{j_l}\beta_le'_l=\beta_l^{-1}a_{j_l}e_l
\end{equation}
so assertion (i) follows for the index $l$, by setting $\phi(l):=j_l$.

Next, set $M_{l,1}:=A_l$ and
$M_{l,2}:=\bigoplus_{i<l} P_i\oplus\bigoplus_{i>l}A_i$.
The foregoing argument provides a subobject
$\alpha:P'_l:=\Img(a'_{\phi(l)}e'_l)\to M_l=M_{l,1}\oplus M_{l,2}$
such that $p'_l\alpha:P'_l\to M_{l,1}$ is an isomorphism. By claim
\ref{cl_immediate}, it follows that the natural map
$\beta'_{l'}:P'_l\oplus M_{l,2}\to M_l$ is an isomorphism. On the
other hand, \eqref{eq_Krull-Remark} implies that $\beta_l$ restricts
to an isomorphism $\gamma_l:P'_l\isom P_l$. We deduce a commutative
diagram
$$
\xymatrix{ P'_l\oplus M_{l,2} \ar[rr]^-{\beta'_{l'}}
\ar[d]_{\gamma_l\oplus\one_{M_{l,2}}} & &
M_l \ar[d]^{\beta_l} \\
M_{l'} \ar[rr]^-{\beta_{l'}} & & A
}$$
whose vertical arrows and top arrow are isomorphisms.
Thus, the bottom arrow is an isomorphism as well, as
required.
\end{pfclaim}

For every $j\in J$, let $e'_j:B_j\to A$ (resp. $p'_j:A\to B_j$)
be the injection (resp. projection) deduced from a given
isomorphism as in (a), and set $a_j:=e'_jp'_j$ for every $j\in J$.
Then $\sum_{j\in J}a_j=\one_A$, so claim \ref{cl_a-plus-b} yields
a mapping $\phi:I\to J$ such that the following holds.
Set $P_j:=\oplus_{i\in\phi^{-1}(j)}\Img(a_je_i)$ for every $j\in J$;
then the natural morphism $\beta:\oplus_{j\in J}P_j\isom A$ is an
isomorphism. However, clearly $\beta(P_j)$ is a subobject of
$\Img(e'_j)$, for every $j\in J$. Hence, $\phi$ must be a surjection,
and $\beta$ restricts to isomorphisms $P_j\isom\Img(e'_j)$ for every
$j\in J$, whence isomorphisms $P_j\isom B_j$, from which (ii) follows
immediately.
\end{proof}

\begin{remark} There are variants of theorem \ref{th_KRS},
that hold under different sets of assumptions. For instance,
in \cite[Ch.I, \S6, Th.1]{Ga} it is stated that part (ii) of
the theorem still holds when $I$ and $J$ are small (not
necessarily finite) sets, provided that $\cC$ admits a small
set of generators and that all small filtered colimits are
representable and exact in $\cC$. (One can moreover show that
the latter result still holds even in the absence of a
small set of generators for $\cC$.)
\end{remark}

\sset\subsubsection{}
Let us now specialize to the case of the category $A\Mod$,
where $A$ is a (commutative) local ring, say with residue
field $k$. Let $M$ be any finitely generated $A$-module;
arguing by induction on the dimension of $M\otimes_Ak$,
one shows easily that $M$ admits a finite decomposition
$M=\oplus_{i=1}^rM_i$, where $M_i$ is indecomposable for
every $i\leq r$.

\begin{lemma}\label{lem_hens-case}
Let $A$ be a henselian local ring, and $M$ an $A$-module
of finite type. Then $M$ is indecomposable if and only
if\/ $\End_A(M)$ is a local ring.
\end{lemma}
\begin{proof} In view of \eqref{subsec_indecomp}, we can
assume that $M$ is indecomposable, and we wish to prove
that $\End_A(M)$ is local. Thus, let $\phi\in\End_A(M)$.
\begin{claim}\label{cl_integral-endo}
The subalgebra $A[\phi]\subset\End_A(M)$ is integral over $A$.
\end{claim}
\begin{pfclaim} This is standard: choose a finite
system of generators $(f_i)_{1\leq i\leq r}$ for $M$,
then we can find a matrix $\ba:=(a_{ij})_{1\leq i,j\leq r}$
of elements of $A$, such that
$\phi(f_i)=\sum_{j=1}^ra_{ij}f_j$ for every $i\leq r$.
The matrix $\ba$ yields an endomorphism $\psi$ of the free
$A$-module $A^{\oplus r}$; on the other hand, let $e_1,\dots,e_r$
be the standard basis of $A^{\oplus r}$ and define an $A$-linear
surjection $\pi:A^{\oplus r}\to M$ by the rule: $e_i\mapsto f_i$
for every $i\leq r$. Then $\phi\circ\pi=\pi\circ\psi$;
by Cayley-Hamilton, $\psi$ is annihilated by the characteristic
polynomial $\chi(T)$ of the matrix $\ba$, whence $\chi(\phi)=0$.
\end{pfclaim}

Since $A$ is henselian, claim \ref{cl_integral-endo} implies
that $A[\phi]$ decomposes as a finite product of local rings.
If there were more than one non-zero factor in this decomposition,
the ring $A[\phi]$ would contain an idempotent $e\neq 1,0$,
whence the decomposition $M=eM\oplus(1-e)M$, where both summands
would be non-zero, which contradicts the assumption.
Thus, $A[\phi]$ is a local ring, so that either $\phi$ or
$\one_M-\phi$ is invertible. Since $\phi$ was chosen arbitrarily
in $\End_A(M)$, the claim follows.
\end{proof}

\begin{corollary}\label{cor_KRS-hens}
Let $A$ be a henselian local ring. Then:
\begin{enumerate}
\item
If $(M_i)_{i\in I}$ and $(N_j)_{j\in J}$ are
two finite families of indecomposable $A$-modules of finite
type such that $\oplus^r_{i=1}M_i\simeq\oplus^s_{j=1}N_j$,
then there is a bijection $\beta:I\isom J$
such that $M_i\simeq N_{\beta(i)}$ for every $i\in I$.
\item
If $M$ and $N$ are two finitely generated $A$-modules
such that $M^{\oplus k}\simeq N^{\oplus k}$ for some
integer $k>0$, then $M\simeq N$.
\item
If $M$, $N$ and $X$ are three finitely generated $A$-modules
such that $X\oplus M\simeq X\oplus N$, then $M\simeq N$.
\end{enumerate}
\end{corollary}
\begin{proof} It follows easily from theorem \ref{th_KRS}
and lemma \ref{lem_hens-case}; the details are left to the
reader.
\end{proof}

\begin{proposition}\label{prop_descend-isom}
Let $A\to B$ be a faithfully flat map of local rings, $M$
and $N$ two finitely presented $A$-modules with a $B$-linear
isomorphism $\omega:B\otimes_AM\isom B\otimes_AN$.
Then $M\simeq N$.
\end{proposition}
\begin{proof} Under the standing assumptions, the natural
$B$-linear map
$$
B\otimes_A\Hom_A(M,N)\to\Hom_B(B\otimes_AM,B\otimes_AN)
$$
is an isomorphism (\cite[Lemma 2.4.29(i.a)]{Ga-Ra}). Hence
we can write $\omega=\sum_{i=1}^r b_i\otimes\phi_i$ for some
$b_i\in B$ and $\phi_i:M\to N$ ($1\leq i\leq r$). Denote by
$k$ and $K$ the residue fields of $A$ and $B$ respectively;
we set
$\bar\phi_i:=\one_k\otimes_A\phi_i:k\otimes_AM\to k\otimes_AN$.
From the existence of $\omega$ we deduce easily that
$n:=\dim_k k\otimes_AM=\dim_k k\otimes_AN$. Hence, after
choosing bases, we can view $\bar\phi_1,\dots,\bar\phi_r$ as
endomorphisms of the $k$-vector space $k^{\oplus n}$.
We consider the polynomial
$p(T_1,\dots,T_r):=\det(\sum^r_{i=1}T_i\cdot\bar\phi_i)\in
k[T_1,\dots,T_r]$. Let $\bar b_1,\dots,\bar b_r$ be the
images of $b_1,\dots,b_r$ in $K$; it follows that
$p(\bar b_1,\dots,\bar b_r)\neq 0$, especially
$p(T_1,\dots,T_r)\neq 0$.

\begin{claim}\label{cl_ok-infinite}
The proposition holds if $k$ is an infinite field or if $k=K$.
\end{claim}
\begin{pfclaim} Indeed, in either of these cases we can
find $\bar a_1,\dots,\bar a_r\in k$ such that
$p(\bar a_1,\dots,\bar a_r)\neq 0$. For every $i\leq r$
choose an arbitrary representative $a_i\in A$ of $\bar a_i$,
and set $\phi:=\sum^r_{i=1}a_i\phi_i$. By construction,
$\one_k\otimes_A\phi:k\otimes_AM\to k\otimes_AN$ is
an isomorphism. By Nakayama's lemma we deduce that
$\phi$ is surjective. Exchanging the roles of $M$
and $N$, the same argument yields an $A$-linear
surjection $\psi:N\to M$. Finally,
\cite[Ch.1, Th.2.4]{Mat} shows that both
$\psi\circ\phi$ and $\phi\circ\psi$ are isomorphisms,
whence the claim.
\end{pfclaim}

Let $A^\he$ be the henselization of $A$, and $A^\sh$, $B^\sh$
the strict henselizations of $A$ and $B$ respectively. Now,
the induced map $A^\sh\to B^\sh$ is faithfully flat, the
residue field of $A^\sh$ is infinite, and $\omega$ induces
a $B^\sh$-linear isomorphism
$B^\sh\otimes_AM\isom B^\sh\otimes_AN$.
Hence, claim \ref{cl_ok-infinite} yields an $A^\sh$-linear
isomorphism
$\beta:A^\sh\otimes_AM\isom A^\sh\otimes_AN$.
However, $A^\sh$ is the colimit of a filtered family
$(A_\lambda~|~\lambda\in\Lambda)$ of finite {\'e}tale
$A^\he$-algebras, so $\beta$ descends to an $A_\lambda$-linear
isomorphism
$\beta_\lambda:A_\lambda\otimes_AM\simeq A_\lambda\otimes_AN$ for
some $\lambda\in\Lambda$. The $A^\he$-module $A_\lambda$
is free, say of finite rank $n$, hence $\beta_\lambda$ can
be regarded as an $A^\he$-linear isomorphism
$(A^\he\otimes_AM)^{\oplus n}\isom(A^\he\otimes_AN)^{\oplus n}$.
Then $A^\he\otimes_AM\simeq A^\he\otimes_AN$, by corollary
\ref{cor_KRS-hens}(ii). Since the residue field of $A^\he$
is $k$, we conclude by another application of claim
\ref{cl_ok-infinite}.
\end{proof}

\sset\subsubsection{Injective hulls}
The notion of injective hull plays a central role in the theory
of local duality : for a noetherian local ring, one constructs
a dualizing module as the injective hull of the residue field
(see \cite[Exp.IV, Th.4.7]{SGA2}), and in section
\ref{sec_homol-polytope}, injective hulls of the residue fields
of a monoid algebra will also enable us to perform a certain
computation of local cohomology, which is a crucial step in the
proof of Hochster's theorem. We present here the basic results
on injective hulls, in the context of arbitrary abelian categories.

\begin{definition}
Let $\cA$ be any abelian category, and $f:N\to M$ a monomorphism
in $\cA$.
\begin{enumerate}
\item
We say that $M$ is an {\em essential extension of $N$} if the
following holds. For any subobject $P\subset M$ we have either
$P=0$ or $P\cap\Img\,f\neq 0$ (here $0$ denotes the zero object
of $\cA$ : see remark \ref{rem_additive-cat}(i)).
\item
We say that $M$ is a {\em proper essential extension of $N$}
if it is an essential extension of $N$, and $f$ is not an
isomorphism.
\item
We say that $M$ is an {\em injective hull\/} of $N$, if $M$
is both an essential extension of $N$, and an injective
object of $\cA$.
\end{enumerate}
\end{definition}

\begin{lemma}\label{lem_hull}
Let $\cA$ be an abelian category, and $I$ an object of $\cA$.
Suppose that :
\begin{enumerate}
\alphaenu
\item
$\cA$ is cocomplete.
\item
All colimits of $\cA$ are universal (see example
{\em\ref{ex_universal-col}}).
\end{enumerate}
Then we have :
\begin{enumerate}
\item
$I$ is injective if and only if it does not admit any proper
essential extensions.
\item
If $N\to M$ is any monomorphism, the set of all essential
extensions of $N$ contained in $M$ admits maximal elements.
\end{enumerate}
\end{lemma}
\begin{proof}(i): Suppose that $I$ is injective, and let $f:I\to M$
be any monomorphism which is not an isomorphism. Then $f$ admits
a left inverse, so $I$ is a direct summand of $M$, hence $M$ is
not an essential extension of $I$.
Conversely, suppose that $I$ does not admit proper essential
extensions; let $f:N\to M$ be a monomorphism in $\cA$ and
$g:N\to I$ any morphism in $\cA$. We consider the cocartesian
diagram in $\cA$
$$
\xymatrix{ N \ar[r]^-f \ar[d]_g & M \ar[d]^{g'} \\
I \ar[r]^-{f'} & P
}$$
and notice that $f'$ is a monomorphism, since the same holds
for $f$. By Zorn's lemma -- and due to conditions (a) and (b) --
we may find a maximal subobject $Q\subset P$ such that
$Q\cap I=0$, and clearly $P/Q$ is an essential extension of $I$
(via $h$); therefore $P/Q=I$, so $P=I\oplus Q$, and if we let
$p:P\to I$ be the resulting projection, we see that
$p\circ g':M\to I$ is an extension of $f$. This shows that $I$
is injective.

(ii): By Zorn's lemma, it suffices to check the following
assertion. Suppose that $(P_j~|~j\in J)$ is a totally ordered
family of essential extensions of $N$ contained in $N$ (for
some totally ordered small indexing set $J$). Then
$P:=\bigcup_{j\in J}P_j$ is again an essential extension of $N$
(notice that this union ({\em i.e.} colimit) is a subobject
of $M$, due to condition (b)). However, say that $Q\subset P$
and $Q\cap N=0$; due to condition (b), we have
$Q=\bigcup_{j\in J}(Q\cap P_j)$, and clearly $(Q\cap P_j)\cap N=0$,
so $Q\cap P_j=0$ for every $j\in J$, and finally $Q=0$.
\end{proof}

\begin{proposition}\label{prop_hull}
Let $\cA$ be an abelian category fulfilling condition {\em (a)}
and {\em (b)} of lemma {\em\ref{lem_hull}}, and $M$ any object
of $\cA$. We suppose moreover that :
\begin{enumerate}
\alphaenu
\addenu\addenu
\item
$\cA$ has enough injective objects.
\end{enumerate}
Then we have :
\begin{enumerate}
\item
$M$ admits an injective hull. More precisely, if $M\to I$ is
any monomorphism into an injective object, then a maximal
essential extension of $M$ in $I$ is an injective hull of $M$.
\item
Let $f:M\to E$ be an injective hull of $M$, and $g:M\to I$
any monomorphism into an injective object. Then $g$ factors
through $f$ and a monomorphism $h:E\to I$.
\item
If $f:M\to E$ and $g:M\to E'$ are two injective hulls of $M$,
there exists an isomorphism $h:E\isom E'$ such that $h\circ f=f'$.
\end{enumerate}
\end{proposition}
\begin{proof}(i): Let $E\subset I$ be such a maximal essential
extension of $M$ (which exists by lemma \ref{lem_hull}(ii)); by
virtue of lemma \ref{lem_hull}(i), it suffices to check that $E$
does not admit proper essential extensions. However, suppose
that $E\to E'$ is a proper essential extension; since $I$ is
injective, the inclusion morphism $E\to I$ extends to a morphism
$f:E'\to I$. By maximality of $E$, we must then have $\Ker\,f\neq 0$;
on the other hand, obviously $E\cap\Ker\,f=0$, which is absurd,
since $E'$ is an essential extension of $E$.

(ii): Since $I$ is injective, $g$ extends to a morphism $h:E\to I$,
and clearly $M\cap\Ker\,h=0$. Since $E$ is an essential extension
of $M$, it follows that $\Ker\,h=0$.

(iii): By (ii) there exists a monomorphism $h:E\to E'$ such
that $h\circ f=f'$. Since $E$ is injective, it follows that
$\Img\,h$ is a direct summand of $E'$; but $E'$ is an essential
extension of $M$, so $\Img\,h=E'$, {\em i.e.} $h$ is also an
epimorphism, hence an isomorphism.
\end{proof}

\sset\subsubsection{}
We specialize now to the case where $\cA$ is the category
$A\Mod$ of $A$-modules, with $A$ an arbitrary noetherian ring.
Recall that if $M$ is any $A$-module, then the set $\Ass\,M$
of all {\em associated primes\/} of $M$ consists of the prime
ideals $\fp\subset A$ such that there exists $m\in M$ with
$\Ann_A(m)=\fp$. (This is the correct definition only for
noetherian rings : we shall see in definition \ref{def_Ass}
a more general notion that is well behaved for arbitrary rings.)
By proposition \ref{prop_hull}(i,iii) the injective hull of $M$
exists and is well defined up to (in general, non-unique)
isomorphism, and we shall denote by $E_A(M)$ a choice of such
hull. 

\begin{lemma}\label{lem_sum-of-inject}
Let $A$ be a noetherian ring. The following holds :
\begin{enumerate}
\item
If $(I_\lambda~|~\lambda\in\Lambda)$ is a (small) family of injective
$A$-modules, then $I:=\bigoplus_{\lambda\in\Lambda}I_\lambda$ is an
injective $A$-module.
\item
If $S\subset A$ is any multiplicative subset, and $M$ any
injective $A$-module, then $M_S$ is an injective $A_S$-module.
\end{enumerate}
\end{lemma}
\begin{proof}(i): Let us first recall :

\begin{claim}\label{cl_prove-later}
Let $R$ be any ring, $M$ an $R$-module. The following
conditions are equivalent :
\begin{enumerate}
\alphaenu
\item
$M$ is an injective $R$-module.
\item
For every ideal $J\subset R$, every $R$-linear map $J\to M$
extends to an $R$-linear map $R\to M$.
\item
$\Ext^1_R(R/J,M)=0$ for every ideal $J\subset R$.
\end{enumerate}
\end{claim}
\begin{pfclaim} Of course, (a)$\Rightarrow$(b)$\Rightarrow$(c).
To check that (c)$\Rightarrow$(a), let $N\subset P$ be an inclusion
of $R$-modules, and $f:N\to M$ an $A$-linear map. We let $S$
be the set of all pairs $(N',f')$, where $N'\subset P$ is an
$R$-submodule containing $N$, and $f':N'\to M$ an $R$-linear
map extending $f$. The set $S$ is partially ordered, by declaring
that $(N',f')\geq(N'',f'')$ if $N''\subset N'$ and $f'_{|N''}=f''$,
for any two pairs $(N',f'),(N'',f'')\in S$. By Zorn's lemma,
$S$ admits a maximal element $(Q,g)$. Suppose $Q\neq P$, and
let $x\in P\setminus Q$. Set $Q':=Q+Rx$ and $J:=\Ann_R(Q'/Q)$;
there follows an exact sequence of $R$-modules
$$
0\to Q\to Q'\to R/J\to 0
$$
and then (c) implies that the restriction map
$\Hom_R(Q',M)\to\Hom_R(Q,M)$ is surjective; especially,
$g$ extetnds to an $R$-linear map $Q'\to M$, contradicting
the maximality of $Q$. Hence $Q=P$, which shows (a).
\end{pfclaim}

In view of claim \ref{cl_prove-later}, it suffices to show that
every $A$-linear map $f:J\to I$ from any ideal $J\subset A$, extends
to an $A$-linear map $A\to I$. However, since $J$ is finitely
generated, there exists a finite subset $\Lambda'\subset\Lambda$
such that the image of $f$ is contained in
$I':=\bigoplus_{\lambda\in\Lambda'}I_\lambda$. Clearly $I'$ is
an injective $A$-module, so $f$ extends to an $A$-linear map
$A\to I'$, and the assertion follows.

(ii): Let $J\subset A_S$ be any ideal, and write $J=I_S$
for some ideal $I\subset A$; since $A$ is noetherian, we have
$$
\Ext^1_{A_S}(A_S/J,M_S)=A_S\otimes_A\Ext^1_A(A/I,M)
$$
so the assertion follows from claim \ref{cl_prove-later}.
\end{proof}

We may now state :

\begin{proposition}\label{prop_decompose-injective}
Let $A$ be any noetherian ring, $M$ any $A$-module. We have :
\begin{enumerate}
\item
For every $\fp\in\Spec\,A$, the $A$-module $E_A(A/\fp)$ is
indecomposable (see \eqref{subsec_indecomp}).
\item
Let $I$ be any non-zero injective $A$-module, and\/ $\fp\in\Ass\,I$
any associated prime. Then $E_A(A/\fp)$ is a direct summand of
$I$. Especially, if $I$ is indecomposable, then $I$ is isomorphic
to $E_A(A/\fp)$.
\item
$\Ass_AM=\Ass_AE_A(M)$.
\item
If\/ $\fp,\fq\in\Spec\,A$ are any two prime ideals, then the
$A$-modules $E_A(A/\fp)$ and $E_A(A/\fq)$ are isomorphic if and
only if\/ $\fp=\fq$.
\item
$E_{A_S}(M_S)\simeq A_S\otimes_AE_A(M)$ for any multiplicative
subset $S\subset A$.
\end{enumerate}
\end{proposition}
\begin{proof}(i): Set $B:=A/\fp$, and suppose that $E_A(B)$ is
decomposable; especially, there exist non-zero submodules
$M_1,M_2\subset E_A(B)$ with $M_1\cap M_2=0$. Thus,
$(M_1\cap B)\cap(M_2\cap B)=0$, and since $B$ is a domain,
we deduce that $M_i\cap B=0$ for $i=1,2$. But since $E_A(B)$
is an essential extension of $B$, this is absurd.

(ii): By assumption, there exists $x\in I$ such that $\Ann_A(x)=\fp$,
so the submodule $Ax\subset I$ is isomorphic to $A/\fp$; then, by
proposition \ref{prop_hull}(ii) there exists a monomorphism
$f:E_A(A/\fp)\to I$, and since $E_A(A/\fp)$ in injective, the image
of $f$ is a direct summand of $I$.

(iii): Clearly $\Ass_A(M)\subset\Ass_AE_A(M)$. Conversely,
suppose $\fp\in\Ass_AE_A(M)$; then there exists an $A$-submodule
$N\subset E_A(M)$ isomorphic to $A/\fp$. Since $E_A(M)$ is an
essential extension of $M$, we have $N\cap M\neq 0$, so
$\fp\in\Ass_A(M)$.

(iv) follows directly from (iii).

(v): We know already that $E':=A_S\otimes_AE_A(M)$ contains $M_S$
and is injective (lemma \ref{lem_sum-of-inject}(iii)), so it
remains only to check that $E'$ is an essential extension of $M$.
Thus, let $x\in E'\setminus M_S$; we have to check that
$N:=A_Sx\cap M_S\neq 0$, and clearly we may assume that $x\in E_A(M)$.
Set $\cF:=\{\Ann_A(tx)~|~t\in S\}$; since $A$ is noetherian,
$\cF$ admits maximal elements, and notice that $A_Sx=A_Stx$
for any $t\in S$. Hence, we may replace $x$ by $tx$ for some
$t\in S$, after which we may assume that $\Ann_A(x)$ is maximal
in $\cF$. We may write $Ax\cap M=Ix$ for some ideal $I\subset A$,
so $N=I_Sx$; let $a_1,\dots,a_k\in A$ be a system of generators
for $I$, and notice that $Ix\neq 0$, since $E_A(M)$ is an essential
extension of $M$. Suppose that $N=0$; then there exists $t\in S$
such that the identity $ta_ix=0$ holds in $A$ for $i=1,\dots,k$.
However, $\Ann_A(x)=\Ann_A(tx)$ by construction, so $Ix=0$, a
contradiction. 
\end{proof}

\begin{theorem}\label{th_classify-injectives}
Let $A$ be a noetherian ring, $I$ an injective $A$-module.
We have :
\begin{enumerate}
\item
$I$ decomposes as a direct sum of the form
\set\begin{equation}\label{eq_decompose-injective}
I\simeq\bigoplus_{\fp\in\Spec\,A}E_A(A/\fp)^{(R_\fp)}
\end{equation}
for a system of (small) sets $(R_\fp~|~\fp\in\Spec\,A)$.
\item
Moreover, the cardinality of $R_\fp$ equals
$\dim_{\kappa(\fp)}\Hom_{A_\fp}(\kappa(\fp),I_\fp)$, for every
$\fp\in\Spec\,A$ (where $\kappa(\fp):=A_\fp/\fp_\fp$). Especially,
this cardinality is independent of the decomposition
\eqref{eq_decompose-injective}.
\end{enumerate}
\end{theorem}
\begin{proof}(i): Denote by $\cF$ the set of all indecomposable
injective submodules of $I$, and by $\cS$ the set of all
subsets $\cG\subset\cF$ such that the natural map
$I_\cG:=\bigoplus_{G\in\cG}G\to I$ is injective. The set $\cS$
is partially ordered by inclusion, and by Zorn's lemma, $\cS$
admits a maximal element $\cM$; in light of proposition
\ref{prop_decompose-injective}(ii), it suffices to check
that $I_\cM=I$. However, $I_\cM$ is injective (lemma
\ref{lem_sum-of-inject}(i)), hence $I=I_\cM\oplus J$ for some
$A$-module $J$, and it is easily seen that $J$ is injective
as well. We are thus reduced to showing that $J=0$. Now, if
$J\neq 0$, let $\fp\in\Ass_AJ$ be any associated prime
(\cite[Th.6.1(i)]{Mat}); then $E_A(A/\fp)$ is an
indecomposable injective direct summand of $J$ (proposition
\ref{prop_decompose-injective}(i,ii)), and therefore
$I_\cM\oplus E_A(A/\fp)$ is a submodule of $I$, contradicting
the maximality of $\cM$.

(ii): In light of (i) and proposition
\ref{prop_decompose-injective}(iii,v) we have
$$
\Hom_{A_\fp}(\kappa(\fp),I_\fp)=
\Hom_{A_\fp}(\kappa(\fp),E_A(A/\fp)^{(R_\fp)}_\fp)=
\kappa(\fp)^{(R_\fp)}\otimes_{\kappa(\fp)}
\Hom_{A_\fp}(\kappa(\fp),E_{A_\fp}(\kappa(\fp)))
$$
so we are reduced to showing the following :

\begin{claim}\label{cl_hom-to-hull}
Let $\fm\subset A$ be any maximal ideal, and set
$\kappa:=A/\fm$; then
$$
d:=\dim_\kappa\Hom_A(\kappa,E_A(\kappa))=1.
$$
\end{claim}
\begin{pfclaim}[] Since $\kappa\subset E_A(\kappa)$, obviously
$d\geq 1$. On the other hand, we have a natural identification
$$
\Hom_A(\kappa,E_A(\kappa))=F:=\{x\in E_A(\kappa)~|~\fm x=0\}.
$$
If $d>1$, we may find $x\in F$ such that $Ax\cap\kappa=0$; but
this is absurd, since $E_A(\kappa)$ is an essential extension
of $\kappa$.
\end{pfclaim}
\end{proof}

\sset\subsubsection{}\label{subsec_coh-inject}
When dealing with the more general coherent rings that appear in
sections \ref{sec_sch-val-rings} and \ref{sec_loc-duality}, the
injective hull is no longer suitable for the study of local
cohomology, but it turns out that one can use instead a
``coh-injective hull'', which works just as well.

\begin{definition}\label{def_coh}
Let $A$ be a ring; we denote by $A\Mod_\coh$ the full subcategory
of $A\Mod$ consisting of all coherent $A$-modules.
\begin{enumerate}
\item
An $A$-module $J$ is said to be {\em coh-injective\/} if the functor
$$
A\Mod_\coh\to A\Mod^o\quad :\quad M\mapsto\Hom_A(M,J)
$$
is exact.
\item
An $A$-module $M$ is said to be {\em $\omega$-coherent\/} if it
is countably generated, and every finitely generated submodule
of $M$ is finitely presented.
\end{enumerate}
\end{definition}

\begin{lemma}\label{lem_omega-coh} Let $A$ be a ring.
\begin{enumerate}
\item
Let $M$, $J$ be two $A$-modules, $N\subset M$ a submodule.
Suppose that $J$ is coh-injective and both $M$ and $M/N$
are $\omega$-coherent. Then the natural map:
$$
\Hom_A(M,J)\to\Hom_A(N,J)
$$
is surjective.
\item
Let $(J_\lambda~|~\lambda\in\Lambda)$ be a filtered family of
coh-injective $A$-modules. Then $\colim_{\lambda\in\Lambda}J_\lambda$
is coh-injective.
\item
Assume that $A$ is coherent, let $M_\bullet$ be an object of\/
$\sD^-(A\Mod_\coh)$ and $I^\bullet$ a bounded below complex
of coh-injective $A$-modules. Then the natural map
$$
\Hom^\bullet_A(M_\bullet,I^\bullet)\to
R\Hom^\bullet_A(M_\bullet,I^\bullet)
$$
is an isomorphism in $\sD(A\Mod)$.
\end{enumerate}
\end{lemma}
\begin{proof}(i): Since $M$ is $\omega$-coherent, we may write it
as an increasing union $\bigcup_{n\in\N}M_n$ of finitely generated,
hence finitely presented submodules. For each $n\in\N$, the image
of $M_n$ in $M/N$ is a finitely generated, hence finitely presented
submodule; therefore $N_n:=N\cap M_n$ is finitely generated,
hence it is a coherent $A$-module, and clearly $N=\bigcup_{n\in\N}N_n$.
Suppose $\phi:N\to J$ is any $A$-linear map; for every $n\in\N$, set
$P_{n+1}:=M_n+N_{n+1}\subset M$, and denote by $\phi_n:N_n\to J$ the
restriction of $\phi$. We wish to extend $\phi$ to a linear map
$\phi':M\to J$, and to this aim we construct  inductively a
compatible system of extensions $\phi'_n:P_n\to J$, for every $n>0$.
For $n=1$, one can choose arbitrarily an extension $\phi'_1$ of
$\phi_1$ to $P_1$. Next, suppose $n>0$, and that $\phi'_n$ is
already given; since $P_n\subset M_n$, we may extend $\phi'_n$
to a map $\phi''_n:M_n\to J$. Since $M_n\cap N_{n+1}=N_n$, and
since the restrictions of $\phi''_n$ and $\phi_{n+1}$ agree on
$N_n$, there exists a unique $A$-linear map $\phi'_{n+1}:P_{n+1}\to J$
that extends both $\phi''_n$ and $\phi_{n+1}$.

(ii): Recall that a coherent $A$-module is finitely presented.
However, it is well known that an $A$-module $M$ is finitely
presented if and only if the functor $Q\mapsto\Hom_A(M,Q)$
on $A$-modules, commutes with filtered colimits (see {\em e.g.}
\cite[Prop.2.3.16(ii)]{Ga-Ra}). The assertion is an easy consequence.

(iii): Since $A$ is coherent, one can find a resolution
$\phi:P_\bullet\to M_\bullet$ consisting of free $A$-modules
of finite rank (cp. \cite[\S7.1.20]{Ga-Ra}). One looks at the
spectral sequences
$$
\begin{aligned}
E_1^{pq}:\Hom_A(P_p,I^q) \Rightarrow\: &
R^{p+q}\Hom^\bullet_A(M_\bullet,I^\bullet) \\
F_1^{pq}:\Hom_A(M_p,I^q) \Rightarrow\: &
H^{p+q}\Hom^\bullet_A(M_\bullet,I^\bullet).
\end{aligned}
$$
The resolution $\phi$ induces a morphism of spectral sequences
$F_1^{\bullet\bullet}\to E_1^{\bullet\bullet}$; on the other
hand, since $I^q$ is coh-injective, we have :
$E_2^{pq}\simeq\Hom_A(H_pF_\bullet,I^q)\simeq
\Hom_A(H_pM_\bullet,I^q)\simeq F^{pq}_2$, so the induced map
$F_2^{pq}\to E_2^{pq}$ is an isomorphism for every $p,q\in\N$,
whence the claim.
\end{proof}

\sset\subsubsection{}\label{eq_SGA2-for-coherents}
Let $A$ be a coherent ring, $Z\subset\Spec\,A$ a constructible
closed subset; we denote by
$$
A\Mod_{\coh,Z}
$$
the full subcategory of $A\Mod_\coh$ whose objects are the
(coherent) $A$-modules with support contained in $Z$.
Let $I\subset A$ be any finitely generated ideal such that
$Z=\Spec\,A/I$; it is easily seen that
$$
\Ob(A\Mod_{\coh,Z})=\bigcup_{n\in\N}\Ob(A/I^n\Mod_\coh).
$$
Now, consider a functor
$$
T:A\Mod_{\coh,Z}^o\to\Z\Mod.
$$
Notice that $TM$ is naturally an $A$-module, for every
coherent $A$-module $M$ : indeed, if $a\in A$ is any
element, we may define the scalar multiplication by $a$
on $TM$ as the $\Z$-linear endomorphism $T(a\cdot\one_M)$.
In other words, $T$ factors through the forgetful functor
$A\Mod\to\Z\Mod$; especially, if we set
$$
H_T:=\colim_{n\in\N}T(A/I^n)
$$
we get a natural $A$-linear map
$$
M\isom\colim_{n\in\N}\Hom_A(A/I^n,M)\to
\colim_{n\in\N}\Hom_A(TM,T(A/I^n))\to
\Hom_A(TM,H_T)
$$
whence a bilinear pairing
$$
M\times TM\to H_T
$$
which in turns yields a natural transformation
$$
\omega_M:TM\to\Hom_A(M,H_T)
\qquad
\text{for every $M\in\Ob(A\Mod_{\coh,Z})$}.
$$

\begin{lemma}\label{lem_SGA2-for-coherents}
In the situation of \eqref{eq_SGA2-for-coherents}, the
following conditions are equivalent :
\begin{enumerate}
\alphaenu
\item
$\omega$ is an isomorphism of functors $T\isom\Hom_A(-,H_T)$.
\item
$T$ is left exact.
\end{enumerate}
\end{lemma}
\begin{proof} Clearly (a)$\Rightarrow$(b). For the converse,
suppose first that $Z=\Spec\,A$, in which case, notice that
$H_T=TA$ and $\omega_A$ is the natural isomorphism
$TA\isom\Hom_A(A,TA)$; then, if $M$ is any coherent $A$-module,
pick a finite presentation
$$
\Sigma
\quad:\quad
A^{\oplus n}\to A^{\oplus m}\to M\to 0
$$
and apply the $5$-lemma to the resulting ladder of $A$-modules
$\omega_\Sigma$ with left exact rows, to deduce that $\omega_M$
is an isomorphism. Next, for a general ideal $I$ and $M$ an
$A$-module with $\Supp\,M\subset Z$, we may find $n\in\N$ such
that $M$ is an $A/I^n$-module; since $A/I^k$ is coherent, the
foregoing case then shows that the induced map
$$
TM\to\Hom_A(M,T(A/I^k))
$$
is an isomorphism, for every $k\geq n$. To conclude, it suffices
to remark that the natural map
$$
\colim_{k\in\N}\Hom_A(M,T(A/I^k))\to\Hom_A(M,H_T)
$$
is an isomorphism, since $M$ is finitely presented
(\cite[Prop.2.3.16(ii)]{Ga-Ra}).
\end{proof}

We wish next to present a criterion that allows to detect,
among the functors $T$ as in \eqref{eq_SGA2-for-coherents},
those that are exact. A complete characterization shall
be given only for a restricted class of coherent ring; namely,
we make the following :

\begin{definition}\label{def_Artin-Rees-rings}
Let $A$ be a coherent ring. We say that $A$ is an
{\em Artin-Rees} ring, if the following holds.
For every finitely generated ideal $I\subset A$, every
coherent $A$-module $M$, and every finitely generated
$A$-submodule $N\subset M$, the $I$-adic topology of $M$
induces the $I$-adic topology on $N$.
\end{definition}

We can then state :

\begin{proposition}\label{prop_SGA2-for-coherents}
In the situation of \eqref{eq_SGA2-for-coherents} suppose
furthermore that $A$ is an Artin-Rees ring. Then the
following conditions are equivalent :
\begin{enumerate}
\alphaenu
\item
$H_T$ is a coh-injective $A$-module.
\item
$T$ is exact.
\end{enumerate}
\end{proposition}
\begin{proof} Clearly (a)$\Rightarrow$(b). For the converse,
let $M$ be any coherent $A$-module, and $N\subset M$ any
finitely generated $A$-submodule; it suffices to check that
the induced map
$$
\Hom_A(M,H_T)\to\Hom_A(N,H_T)
$$
is surjective. However, let $f:N\to H_T$ be any $A$-linear map;
since $N$ is finitely presented, there exists $n\in\N$ such that
$f$ factors through an $A$-linear map $f_n:N\to T(A/I^n)$
(\cite[Prop.2.3.16(ii)]{Ga-Ra}), and clearly $I^nN\subset\Ker\,f_n$,
so $I^nN\subset\Ker\,f$ as well.
Since the $I$-adic topology of $N$ agrees with the topology
induced by the $I$-adic topology of $M$, there exists $k\in\N$
such that $I^kM\cap N\subset I^nN$.
Set $\bar N:=N/(I^kM\cap N)$; then $\bar N$ is a finitely
generated $A$-submodule of $M/I^kM$, and especially it is
a coherent $A$-module. By construction, $f$ factors through
an $A$-linear map $\bar f:\bar N\to H_T$; assumption (b)
and lemma \ref{lem_SGA2-for-coherents} imply that $\bar f$
extends to an $A$-linear map $M/I^kM\to H_T$, and the
resulting map $M\to H_T$ extends $f$, whence (a).
\end{proof}

\begin{remark}\label{rem_SGA2-for-coherents}
(i)\ \
In the terminology of definition \ref{def_Artin-Rees-rings},
the standard Artin-Rees lemma implies that every noetherian
ring is an Artin-Rees ring. We shall see later that if $V$
is any valuation ring, then every essentially finitely
presented $V$-algebra is an Artin-Rees ring (corollary
\ref{cor_coherence} and theorem \ref{th_Rees}).

(ii)\ \
On the other hand, if $A$ is noetherian, claim \ref{cl_prove-later}
easily implies that an $A$-module is coh-injective if and only
if it is injective.

(iii)\ \
Combining (i) and (ii) with proposition \ref{prop_SGA2-for-coherents},
we recover \cite[Exp.IV, Prop.2.1]{SGA2}.
\end{remark}

\begin{example}\label{ex_SGA2-for-coherents}
Let $\kappa_0$ be a field, $A$ a noetherian $\kappa_0$-algebra,
$\fm\subset A$ a maximal ideal such that $\kappa:=A/\fm$ is a
finite extension of $\kappa_0$, and set $Z:=\{\fm\}\subset\Spec\,A$.
Then every object of $A\Mod_{\coh,Z}$ is a finite dimensional
$\kappa_0$-vector space, and therefore the functor
$$
T:A\Mod_{\coh,Z}\to\kappa_0\Mod
\qquad
M\mapsto\Hom_{\kappa_0}(M,\kappa_0)
$$
is exact. Proposition \ref{prop_SGA2-for-coherents} and
remark \ref{rem_SGA2-for-coherents}(ii) then say that
$$
H_T:=\colim_{n\in\N}\Hom_{\kappa_0}(A/\fm^n,\kappa_0)
$$
is an injective $A$-module. More precisely, notice that
$\Hom_A(\kappa,H_T)=\Ann_{H_T}(\fm)=T(\kappa)\simeq\kappa$,
therefore $H_T$ is the injective hull of the residue field
$\kappa$.
\end{example}

\sset\subsubsection{Flatness criteria}
The following generalization of the local flatness
criterion answers affirmatively a question raised in
\cite[Ch.IV, Rem.11.3.12]{EGAIV-3}.

\begin{lemma}\label{lem_Tor-crit}
Let $A$ be a ring, $I\subset A$ an ideal, $B$ a finitely
presented $A$-algebra, $\fp\subset B$ a prime ideal
containing $IB$, and $M$ a finitely presented $B$-module.
Then the following two conditions are equivalent:
\begin{enumerate}
\alphaenu
\item
$M_\fp$ is a flat $A$-module.
\item
$M_\fp/IM_\fp$ is a flat $A/I$-module and $\Tor_1^A(M_\fp,A/I)=0$.
\end{enumerate}
\end{lemma}
\begin{proof} Clearly, it suffices to show that (b)$\Rightarrow$(a),
hence we assume that (b) holds. We write $A$ as the union of
a filtered family $(A_\lambda~|~\lambda\in\Lambda)$ of its
$\Z$-subalgebras of finite type, and set $A':=A/I$,
$I_\lambda:=I\cap A_\lambda$, $A'_\lambda:=A_\lambda/I_\lambda$
for every $\lambda\in\Lambda$.
Then, for some $\lambda\in\Lambda$, the $A$-algebra
$B$ descends to an $A_\lambda$-algebra $B_\lambda$ of finite
type, and $M$ descends to a finitely presented
$B_\lambda$-module $M_\lambda$. For every $\mu\geq\lambda$
we set $B_\mu:=A_\mu\otimes_{A_\lambda}B_\lambda$ and
$M_\mu:=B_\mu\otimes_{B_\lambda}M_\lambda$. Up to replacing
$\Lambda$ by a cofinal family, we can assume that $B_\mu$
and $M_\mu$ are defined for every $\mu\in\Lambda$. Then,
for every $\lambda\in\Lambda$ let
$g_\lambda:B_\lambda\to B$ be the natural map, and set
$\fp_\lambda:=g_\lambda^{-1}\fp$.

\begin{claim}\label{cl_flat-is-ok}
There exists $\lambda\in\Lambda$ such that
$M_{\lambda,\fp_\lambda}/I_\lambda M_{\lambda,\fp_\lambda}$
is a flat $A'_\lambda$-module.
\end{claim}
\begin{pfclaim} Set
$B'_\lambda:=B_\lambda\otimes_{A_\lambda}A'_\lambda$ and
$M'_\lambda:=M_\lambda\otimes_{A_\lambda}A'_\lambda$
for every $\lambda\in\Lambda$; clearly the natural maps
$$
\colim_{\lambda\in\Lambda}A'_\lambda\to A/IA
\qquad
\colim_{\lambda\in\Lambda}B'_\lambda\to B/IB
\qquad
\colim_{\lambda\in\Lambda}M'_\lambda\to M/IM
$$
are isomorphisms. Then the claim follows from (b) and
\cite[Ch.IV, Cor.11.2.6.1(i)]{EGAIV-3}.
\end{pfclaim}

In view of claim \ref{cl_flat-is-ok}, we can replace $\Lambda$
by a cofinal subset, and thereby assume that
$M_{\lambda,\fp_\lambda}/I_\lambda M_{\lambda,\fp_\lambda}$ is a flat
$A'_\lambda$-module for every $\lambda\in\Lambda$.

\begin{claim}\label{cl_colim&Tors}
(i)\ \  
More generally, let $R_\bullet:=(R_\lambda~|~\lambda\in\Lambda)$
be any system of rings indexed by any filtered set $\Lambda$,
and $N_\bullet:=(N_\lambda~|~\lambda\in\Lambda)$,
$N'_\bullet:=(N'_\lambda~|~\lambda\in\Lambda)$ two systems
consisting of $R_\lambda$-modules $N_\lambda$ and $N'_\lambda$
for every $\lambda\in\Lambda$, and $R_\lambda$-linear transition
maps $N_\lambda\to N_\mu$, $N'_\lambda\to N'_\mu$ for every
$\lambda,\mu\in\Lambda$ with $\lambda\leq\mu$. Denote by
$R$, $N$, $N'$ the colimits of $R_\bullet$, $N_\bullet$ and
respectively $N'_\bullet$; then the natural maps
$R_\lambda\to R$, $N_\lambda\to N$, $N'_\lambda\to N$ induce
isomorphisms :
$$
\colim_{\lambda\in\Lambda}\Tor^{R_\lambda}_i(N_\lambda,N'_\lambda)\to
\Tor_i^R(N,N')
\qquad
\text{for every $i\in\N$}.
$$
\begin{enumerate}
\addenu
\item
For every $\lambda,\mu\in\Lambda$ with $\mu\geq\lambda$, the
natural map:
$$
f_{\lambda\mu}:
B_{\mu,\fp_\mu}\otimes_{B_\lambda}\Tor^{A_\lambda}_1(M_{\lambda,\fp_\lambda},A'_\lambda)
\to\Tor^{A_\mu}_1(M_{\mu,\fp_\mu},A'_\mu)
$$
is surjective.
\end{enumerate}
\end{claim}
\begin{pfclaim}(i): For every $\lambda\in\Lambda$, let
$L_\bullet(N'_\lambda)$ be the canonical free resolution of the
$R_\lambda$-module $N'_\lambda$ (\cite[Ch.X, \S3, n.3]{BouAH}).
Similarly, denote by $L_\bullet(N')$ the canonical free
resolution of the $R$-module $N'$. It follows from
\cite[Ch.II, \S6, n.6, Cor.]{Bourbaki} and the exactness
properties of filtered colimits, that the natural map:
$\colim_{\lambda\in\Lambda}L_\bullet(N'_\lambda)\to
L_\bullet(N')$ is an isomorphism. Hence :
$$
\begin{aligned}
\colim_{\lambda\in\Lambda}
H_\bullet(N_\lambda\derotimes_{R_\lambda}N'_\lambda) \simeq\: &
\colim_{\lambda\in\Lambda}
H_\bullet(N_\lambda\otimes_{R_\lambda}L_\bullet(N'_\lambda)) \\
\simeq\: &
H_\bullet(N\otimes_R\colim_{\lambda\in\Lambda}L_\bullet(N'_\lambda)) \\
\simeq\: & H_\bullet(N\otimes_RL_\bullet(N')) \\
\simeq\: & H_\bullet(N\derotimes_RN').
\end{aligned}
$$

(ii): We use the change of ring spectral sequences for $\Tor$
(\cite[Th.5.6.6]{We})
$$
\begin{aligned}
E^2_{pq} :\: &
\Tor^{A_\mu}_p(\Tor^{A_\lambda}_q(M_{\lambda,\fp_\lambda},A_\mu),A'_\mu)
\Rightarrow \Tor^{A_\lambda}_{p+q}(M_{\lambda,\fp_\lambda},A'_\mu) \\
F^2_{pq} :\: &
\Tor^{A'_\lambda}_p(\Tor^{A_\lambda}_q(M_{\lambda,\fp_\lambda},A'_\lambda),A'_\mu)
\Rightarrow \Tor^{A_\lambda}_{p+q}(M_{\lambda,\fp_\lambda},A'_\mu).
\end{aligned}
$$
Since $F^2_{10}=F^2_{20}=0$, the natural map
$$
\alpha:F_{01}^2:=
A'_\mu\otimes_{A'_\lambda}\Tor^{A_\lambda}_1(M_{\lambda,\fp_\lambda},A'_\lambda)
\to\Tor^{A_\lambda}_1(M_{\lambda,\fp_\lambda},A'_\mu)
$$
is an isomorphism. On the other hand, we have a surjection:
$$
\beta:\Tor_1^{A_\lambda}(M_{\lambda,\fp_\lambda},A'_\mu)
\to E^2_{10}:=\Tor^{A_\mu}_1(M_{\mu,\fp_\lambda},A'_\mu)
$$
and $f_{\lambda\mu}=(\beta\circ\alpha)_{\fp_\mu}$, whence the claim.
\end{pfclaim}

We deduce from claim \ref{cl_colim&Tors}(i) that the natural map
$$
\colim_{\lambda\in\Lambda}
\Tor^{A_\lambda}_1(M_{\lambda,\fp_\lambda},A'_\lambda)
\to\Tor_1^A(M_\fp,A')=0
$$
is an isomorphism. However,
$\Tor^{A_\lambda}_1(M_{\lambda,\fp_\lambda},A'_\lambda)$
is a finitely generated $B_{\lambda,\fp_\lambda}$-module by
\cite[Ch.X, \S6, n.4, Cor.]{BouAH}. We deduce that
$f_{\lambda\mu,\fp_\lambda}=0$ for some $\mu\geq\lambda$;
therefore $\Tor^{A_\mu}_1(M_{\mu,\fp_\mu},A'_\mu)=0$, in view
of claim \ref{cl_colim&Tors}(ii), and then the local
flatness criterion of \cite[Ch.0, \S10.2.2]{EGAIII} says that
$M_{\mu,\fp_\mu}$ is a flat $A_\mu$-module, so finally $M_\fp$
is a flat $A$-module, as stated.
\end{proof}

\section{Complements of topology and topological algebra}
\label{chap_top-algebra}
This chapter is a miscellanea of results of topology and
topological algebra that shall be needed in the rest of
the treatise.

\subsection{Spectral spaces and constructible subsets}
The class of spectral topological spaces was first
introduced in Hochster's paper \cite{Hoch2},
where it was shown to coincide precisely with
the class of those spaces that are homeomorphic to
(Zariski) spectra of commutative rings. Since then,
spectral spaces have reappeared in a variety of
unrelated contexts; their ubiquity has all to do with
the very agreeable properties of their boolean algebras
of constructible subsets; especially, the constructible
subsets of a spectral space generate a quasi-compact
and separated topology. The many corollaries springing
from this basic observation add up to a handy toolkit
for dealing in a uniform way with the kind of
non-separated spaces that most frequently occur in
geometric applications of an algebraic bent.

In this section we introduce spectral spaces and
prove their basic properties. We also collect a
few results showing the interplay between algebraic
and topological properties of schemes and their
constructible subsets.

\begin{definition} Let $X$ be a topological space.
\begin{enumerate}
\item
We say that $X$ is {\em quasi-separated} if the
intersection of any two quasi-compact open subsets
of $X$ is still quasi-compact.
\item
We say that $X$ is {\em coherent} if it is quasi-compact
and it admits a basis consisting of quasi-compact open
subsets, and closed under finite intersections.
\item
We say that $X$ is {\em reducible} if there exist non-empty
closed subsets $Z,Z'\subset X$ such that $Z\cup Z'=X$ and
$Z,Z'\neq X$. We say that $X$ is {\em irreducible}
if it is non-empty and not reducible.
\item
We say that $X$ is {\em sober} if for every irreducible
closed subset $Z\subset X$ (with the topology induced
from $X$) there exists a unique point $\eta_Z\in Z$ such
that $Z$ is the topological closure of $\{\eta_Z\}$ in $X$.
In this case, we say that $\eta_Z$ is the
{\em generic point} of $Z$.
\item
We say that $X$ is {\em spectral} if it is sober and
coherent.
\item
We say that $X$ is {\em noetherian}, if for every
descending chain
$$
\cdots\subset Z_2\subset Z_1\subset Z_0
$$
of closed subsets of $X$, there exists $n\in\N$ such that
$Z_m=Z_n$ for every $m\geq n$.
\item
We say that $X$ is {\em locally coherent} (resp.
{\em locally spectral}, resp. {\em locally noetherian})
if every point of $X$ admits an open neighborhood which
is coherent (resp. spectral, resp. noetherian).
\item
We say that $X$ is $T_0$ if, for every pair of points $x,y\in X$,
there exists an open subset $U\subset X$ such that $U\cap\{x,y\}$
contains exactly one point.
\item
We say that $X$ is {\em compact} if it is quasi-compact
and separated.
\end{enumerate}
\end{definition}

\begin{remark}\label{rem_irred-comps}
(i)\ \
Let $X$ be a topological space, $(Z_i~|~i\in I)$ a
family of irreducible closed subsets of $X$, totally
ordered by inclusion, and denote by $Z$ the topological
closure in $X$ of $\bigcup_{i\in I}Z_i$. Then $Z$ is also
irreducible. Indeed, say that $Z=Y_1\cup Y_2$ for two
closed subsets $Y_1,Y_2$ of $X$. If
$Z_i\subset Y_1\cap Y_2$ for every $i\in I$, then
$Z\subset Y_1$; otherwise, say that $Z_i\not\subset Y_1$
for some $i\in I$; then $Z_i\subset Y_2$ for every
$i\in I$, whence $Z\subset Z_2$.

(ii)\ \
By Zorn's lemma, it follows that every irreducible
closed subset of $X$ is contained in a maximal
one. A maximal irreducible closed subset of $X$
is called an {\em irreducible component} of $X$,
and clearly $X$ is the union of its irreducible components.

(iii)\ \
Let $Z_1,Z_2\subset X$ be any two closed subsets,
and endow $Z_1$, $Z_2$ and $Z_3:=Z_1\cup Z_2$ with the
topology induced by $X$. For $j\leq 3$, denote by
$I_j$ the set of all irreducible components of $Z_i$;
then $I_3\subset I_1\cup I_2$. Indeed, let $Y\in I_3$;
then $Y=(Y\cap Z_1)\cup(Y\cap Z_2)$, so we must have
$Y=Y\cap Z_i$ for some $i\in\{1,2\}$.

(iv)\ \
Let $k\in\N$ be any integer, and
$Z_\bullet:=(Z_i~|~i=0,\dots,k)$ a strictly descending
finite chain of irreducible closed subsets of $X$.
The integer $k$ is the {\em length} of $Z_\bullet$.
The supremum of the lengths of all such chains is
an invariant called the {\em (Krull) dimension} of
$X$, and denoted
$$
\dim X
$$
(especially, $\dim\emptyset=-\infty$). We say that $X$ has
{\em finite Krull dimension} if $\dim X\in\N\cup\{-\infty\}$.

(v)\ \
Suppose that $X$ is the union of finitely many closed
subsets $Z_1,\dots,Z_k$. Then it follows directly from
(iii) that $\dim X=\sup(\dim Z_1,\dots,\dim Z_k)$.

(vi)\ \
If $X$ is sober, then $X$ is also $T_0$. Indeed, let
$x,y\in X$ be any two distinct points, and denote by
$Z_y$ the topological closure of $\{y\}$ in $X$; if
$X$ is sober, we may assume -- up to swapping the roles
of $x$ and $y$ -- that $x\notin Z_y$. In this case, the
set $X\setminus Z_y$ is an open neighborhood of $x$ that
does not contain $y$, whence the claim.
\end{remark}

\begin{lemma}\label{lem_irreducible}
Let $X$ be a topological space, $T\subset X$ any subset,
$\bar T$ the topological closure of $T$ in $X$, and endow
$T$ and $\bar T$ with the topologies induced from $X$. The
following holds :
\begin{enumerate}
\item
$T$ is irreducible if and only if the same holds for $\bar T$.
\item
If\/ $T$ is non-empty and open in $X$, and\/ $X$ is irreducible,
then $T$ is irreducible.
\item
If $X$ is sober, and $T$ is locally closed in $X$,
then $T$ is sober as well.
\item
Suppose that $T=T_1\cup T_2$ for two subsets $T_1$
and $T_2$, and endow also $T_1$ and $T_2$ with the
topologies induced from $X$. If\/ $T_1$, $T_2$ and
$X$ are sober, the same holds for $T$.
\item
Let $f:X\to Y$ be any continuous map, and endow $f(X)$
with the topology induced by $Y$. If\/ $X$ is irreducible,
the same holds for $f(X)$.
\end{enumerate}
\end{lemma}
\begin{proof} For any subset $Y\subset X$, let $\bar Y$
denote the topological closure of $Y$ in $X$. 

(i): Suppose that $T$ is irreducible, and say that
$\bar T=V\cup V'$ for two closed subsets $V,V'$ of $\bar T$;
then $T=(T\cap V)\cup(T\cap V')$, so we may assume that
$T=T\cap V$, in which case $\bar T=V$. Conversely, suppose
that $\bar T$ is irreducible, and $T=Z_1\cup Z_2$ for some
closed subsets $Z_1$, $Z_2$ of $T$; then
$\bar T=\bar Z_1\cup\bar Z_2$, so we may assume that
$\bar T=\bar Z_1$, in which case $T=T\cap\bar Z_1=Z_1$.

(ii): Indeed, say that $T=Z_1\cup Z_2$ for two non-empty
closed subsets of $T$; then
$X=\bar Z_1\cup\bar Z_2\cup(X\setminus T)$. By assumption,
$X\not\subset X\setminus T$; as $X$ is irreducible, we may
then assume that $X=\bar Z_1$, whence $T=Z_1$.

(iii): This is clear if $T$ is closed in $X$, so
we may assume that $T$ is open in $X$. However,
say that $Z\subset T$ is an irreducible closed
subset; by (i), the topological closure $\bar Z$
of $Z$ in $X$ is irreducible, and clearly its
generic point  $\eta_{\bar Z}$ lies in $Z$, so
$\eta_{\bar Z}$ is the unique point of $Z$ such
that the closure of $\{\eta_{\bar Z}\}$ in $T$
equals $Z$.

(iv): Let $Z$ be any irreducible closed subset
of $T$, and set $Z_i:=Z\cap T_i$ for $i=1,2$.
According to (i), $\bar Z$ is irreducible; since
$\bar Z=\bar Z_1\cup\bar Z_2$, we may then assume
that $\bar Z=\bar Z_1$, especially $\bar Z_1$ is
irreducible, so $Z_1$ is irreducible, by (i). Therefore
$Z_1$ admits a generic point $\eta$, since $T_1$ is
sober. By construction, the topological closure
of $\{\eta\}$ in $T$ equals $T\cap\bar Z=Z$,
and it remains to see that $\eta$ is the unique
point of $Z$ with this property. However, if
$\eta'$ is any other such point, then clearly
the topological closures in $X$ of both $\{\eta\}$
and $\{\eta'\}$ equal $\bar Z$, so $\eta=\eta'$,
since $X$ is sober.

(v): Indeed, suppose that $f(X)=Z_1\cup Z_2$ for two
closed subsets $Z_1,Z_2$ of $f(X)$; then
$X=f^{-1}(Z_1)\cup f^{-1}(Z_2)$, so we may assume that
$X=f^{-1}Z_1$, in which case $f(X)=Z_1$.
\end{proof}

\begin{proposition}\label{prop_noetherian}
Let $X$ be any topological space. Then $X$ is
noetherian if and only if the following conditions
hold :
\begin{enumerate}
\alphaenu
\item
Every closed subset of $X$ has a finite number of
irreducible components.
\item
For every descending chain $(Z_n~|~n\in\N)$
of irreducible closed subsets of $X$, there exists
$N\in\N$ such that $Z_n=Z_N$ for every $n\geq N$.
\end{enumerate}
\end{proposition}
\begin{proof} Suppose first that $X$ is noetherian.
Then (b) clearly holds. To show (a), let $S$ denote
the set of all closed subsets of $X$ that have
infinitely many irreducible components; since $X$
is noetherian, if $S$ is not empty, it admits
minimal elements, so say that $Z$ is such a minimal
element. Then clearly $Z$ cannot be irreducible, so
$Z=Z_1\cup Z_2$ for some non-empty closed subsets of
$X$ strictly contained in $Z$. By minimality of $Z$,
we have $Z_1,Z_2\notin S$; but then remark
\ref{rem_irred-comps}(iii) implies that
$Z_1\cup Z_2\notin S$ as well, a contradiction.

For the converse, we show more precisely the following :

\begin{claim} Let $X$ be any topological space,
$Z_\bullet:=(Z_n~|~n\in\N)$ a descending non-stationary
chain of closed subsets of $X$, such that $Z_n$ has a
finite number of irreducible components, for every
$n\in\N$. Then there exists a descending, non-stationary
chain $(Z'_n~|~n\in\N)$ of closed subsets of $X$,
such that $Z'_n$ is an irreducible component of
$Z_n$, for every $n\in\N$.
\end{claim}
\begin{pfclaim}[] For every $n\in\N$, let $I_n$ be
the (finite) set of irreducible components of $Z_n$,
and denote by $I'_n\subset I_n$ the subset of all
$T\in I_n$ such that $T\in I_m$ for every $m\geq n$.
Notice that, since $Z_\bullet$ is descending and
non-stationary, the subset $I''_n:=I_n\setminus I'_n$
is non-empty, for every $n\in\N$. Moreover, for every
$n>0$, every $T\in I''_n$ and every $T'\in I_{n-1}$
such that $T\subset T'$, we must have $T'\in I''_{n-1}$;
indeed, otherwise we would have $T'\in I_n$, hence
$T=T'$, hence $T\in I'_n$, a contradiction.
Hence, for every $n\in\N$, denote by $S_n$ the set of
all chains $(T_k~|~k=0,\dots,n)$ such that
\begin{itemize}
\item
$T_k\in I''_k$ for every $k=0,\dots,n$.
\item
$T_k\subset T_{k-1}$ for every $k=1,\dots,n$.
\end{itemize}
The foregoing observations easily imply that
$S_n$ is a finite non-empty set for every $n\in\N$.
Moreover, we have an obvious map $S_m\to S_n$ whenever
$m\geq n$, that assigns to any $T_\bullet\in S_m$
the truncated chain $(T_k~|~k=0,\dots,n)$. The
limit $S$ of the system $(S_n~|~n\in\N)$ is then
non-empty; say that $(T_\bullet^k~|~k\in\N)$ is
any element of $S$, and set $Z'_n:=T^n_n$ for every
$n\in\N$. We claim that the chain $(Z'_n~|~n\in\N)$
will do. Indeed, suppose by way of contradiction,
that the chain $Z'_\bullet$ is stationary, so there
exists $N\in\N$ such that $Z'_n=Z'_N$ for every $n\geq N$;
but then we would have $Z_N\in I_n$ for every $n\geq N$,
hence $Z_N\in I'_N$, which is absurd.
\end{pfclaim}
\end{proof}

\begin{remark}\label{rem_sorite-spectral}
Let $X$ be a topological space, $U\subset X$ an open subset,
and endow $U$ with the topology induced from $X$.

(i)\ \
If $X$ is quasi-separated, then the same holds for $U$.
If $X$ is coherent, then $X$ is quasi-separated. Indeed,
let $(U_i~|~i\in I)$ be a basis for $X$ consisting of
quasi-compact open subsets and closed under finite
intersections, and $U_1$, $U_2$ two quasi-compact open
subsets. Then there exist finite subsets $J_i\subset I$
such that $U_i=\bigcup_{i\in J_i}U_i$ for $i=1,2$; consequently
$U_1\cap U_2=\bigcup_{(j,j')\in J_1\times J_2}(U_j\cap U_{j'})$,
and by assumption each $U_j\cap U_{j'}$ is quasi-compact,
whence the claim.

(ii)\ \
If $X$ is coherent (resp. spectral) and $U$ is
quasi-compact, then $U$ is coherent (resp. spectral)
for the topology induced from $X$. Indeed, lemma
\ref{lem_irreducible}(iii) says that $U$ is sober
whenever the same holds for $X$, and if
$(U_i~|~i\in I)$ is a basis of quasi-compact open
subsets of $X$, then $(U\cap U_i~|~i\in I)$ is a
basis of quasi-compact open subsets of $U$.

(iii)\ \
It follows easily from (ii) that if $X$ is locally
coherent (resp. locally spectral) then the same holds
for $U$.

(iv)\ \
If $X$ is noetherian, every open subset of $X$ is
quasi-compact, so $X$ is coherent. Also, every locally
closed subset of $X$ is noetherian as well.

(v)\ \
Furthermore, suppose that $T_1,\dots,T_k$ is any
finite family of subsets of $X$, such that
$X=\bigcup_{i=1}^kT_i$, and endow each $T_i$ with
the topology induced from $X$; then, if each $T_i$
is noetherian, the same holds for $X$.
\end{remark}

\sset\subsubsection{}
Let us denote by $\Top$ the category of topological spaces
(whose morphisms are the continuous maps), and by
$\mathbf{Sober}$ the full subcategory of $\Top$ whose
objects are the sober topological spaces.

\begin{proposition}\label{prop_soberify}
The inclusion functor
$$
\mathbf{Sober}\to\Top
$$
admits a left adjoint.
\end{proposition}
\begin{proof} For any topological space $X$, denote
by $S_X$ the set of all irreducible closed subsets
of $X$. If $Z$ is any irreducible closed subset of
$X$, we shall use the notation $\eta_Z$ to refer
to $Z$, when we wish to regard the latter as an
element of $S_X$. From lemma \ref{lem_irreducible}(i,ii)
we see that, for every open subset $U\subset X$, we have
a natural bijection
$$
S_U\isom\{\eta_Z\in S_X~|~Z\cap U\neq\emptyset\}
\qquad
T\mapsto\eta_{\bar T}
$$
(where, for any subset $T$ of $X$, we let $\bar T$ be
the topological closure of $T$ in $X$). We may then
identify $S_U$ with its image in $S_X$, and we notice

\begin{claim}\label{cl_endow}
(i)\ \
$S_U\cap S_V=S_{U\cap V}$ for every open subsets
$U,V\subset X$.
\begin{enumerate}
\addenu
\item
For any family $(U_i~|~i\in I)$ of open subsets
of $X$, we have $S_{\bigcup_{i\in I}U_i}=\bigcup_{i\in I}S_{U_i}$.
\end{enumerate}
\end{claim}
\begin{pfclaim}(i): Clearly $S_{U\cap V}\subset S_U\cap S_V$.
For the converse inclusion, suppose that $Z$ is closed
subset of $X$, such that $Z\cap U,Z\cap V\neq\emptyset$
and $Z\cap U\cap V=\emptyset$. Then
$Z=(Z\!\setminus\!U)\cup(Z\!\setminus\!V)$, so $Z$ is
not irreducible in $X$. The assertion is an immediate
consequence.

(ii) is trivial.
\end{pfclaim}

From claim \ref{cl_endow} we see that the family
$(S_U~|~U\ \text{is open in $X$})$ is a topology
on $S_X$, and we endow $S_X$ with this topology.
Next, there is a natural map
$$
f_X:X\to S_X
\qquad
x\mapsto\eta_{\bar{\{x\}}}
$$
and it is easily seen that $f^{-1}_X(S_U)=U$ for
every open subset $U\subset X$; especially $f_X$
is continuous. It also follows that the rule
$Z\mapsto f_X^{-1}Z$ establishes a bijection
from the closed subsets of $S_X$ to those of $X$;
since this bijection respects inclusion of closed
subsets, we deduce that it restricts to a bijection
from the irreducible closed subsets of $S_X$ to
those of $X$. Moreover, in light of lemma
\ref{lem_irreducible}(v), it is easily seen that
the inverse of this bijection assigns to any irreducible
closed subset $Z$ of $X$, the topological closure
of $f_X(Z)$ in $S_X$.

\begin{claim}\label{cl_S-is-sober}
$S_X$ is a sober topological space.
\end{claim}
\begin{pfclaim} Let $Z$ be any irreducible closed
subset of $X$, and denote by $S_Z$ the topological
closure of $\{\eta_Z\}$ in $S_X$; then
$S_Z=S_X\!\setminus\!S_U$, where $S_U$ is the largest
open subset of $S_X$ that does not contain $\eta_Z$,
{\em i.e.} $U$ is the largest open subset of $X$
that does not intersect $Z$. So $U=X\!\setminus\!Z$,
and therefore $f^{-1}_X(S_Z)=Z$. We conclude that
the topological closure of $f_X(Z)$ in $S_X$ equals
$S_Z$. Summing up, we see that the rule
$\eta_Z\mapsto S_Z$ establishes a bijection from
$S_X$ to the set of irreducible closed subsets of
$S_X$, whence the claim.
\end{pfclaim}

Now, let $Y$ be any sober space, and $g:X\to Y$ any
continuous map; in light of claim \ref{cl_S-is-sober},
it suffices to show that there exists a unique continuous
map $h:S_X\to Y$ such that $g=h\circ f_X$. However,
it is easily seen that any continuous map $h$ with
this property must assign to every $\eta_Z\in S_X$
the generic point $\eta_{g(Z)}$ of the topological
closure of $g(Z)$ in $S$ (lemma \ref{lem_irreducible}(v));
conversely, one checks easily that the rule
$$
\eta_Z\mapsto\eta_{g(Z)}
\qquad
\text{for every irreducible closed subset $Z$ of $X$}
$$
defines a map $h:S_X\to Y$ such that $h^{-1}U=S_{g^{-1}U}$
for every open subset $U\subset Y$; especially, $h$
is continuous : the details shall be left to the reader.
\end{proof}

\begin{remark}\label{rem_sober-permanence}
Keep the notation of the proof of proposition
\ref{prop_soberify}, and let $U\subset X$ be any open
subset; it is easily seen that $U=\emptyset$ if and
only if $S_U=\emptyset$. It follows that the rule
$V\mapsto f_X^{-1}V$ defines a bijection from the set
of open subsets of $S_X$ to the set of open subsets
of $X$. Consequently, $X$ is quasi-compact (resp.
quasi-separated) if and only if the same holds for
$S_X$, and $X$ is coherent if and only if $S_X$ is
spectral.
\end{remark}

\sset\subsubsection{}\label{subsec_sober-lim}
Let $X_\bullet:=(X_i~|~i\in\Ob(I))$ be a system of sober
and quasi-compact topological spaces, indexed by a
category $I$, with continuous transition maps
$p_u:X_i\to X_j$ for every morphism $u:i\to j$ in $I$,
and set
$$
X:=\lim_{i\in I}X_i
$$
where the limit is taken in the category of topological
spaces. For every $i\in\Ob(I)$, denote by $p_i:X\to X_i$
the induced projection, and recall that the system of all
subsets of the form $\bigcap_{j\in J}p_j^{-1}U_j$, where $J$
ranges over all finite subsets of $\Ob(I)$, and $U_j$
ranges over all open subsets of $X_j$, is a basis of the
topology of $X$.

\begin{proposition}\label{prop_lim-qc-sober}
In the situation of \eqref{subsec_sober-lim}, the
following holds :
\begin{enumerate}
\item
$X$ is a sober topological space.
\item
Suppose furthermore, that the system $X_\bullet$
is cofiltered. We have :
\begin{enumerate}
\item
If $X_i\neq\emptyset$ for every $i\in\Ob(I)$, then
$X$ is quasi-compact and non-empty.
\item
If $p_u$ is quasi-compact for every morphism $u$
of $I$, then $p_i$ is quasi-compact, for every
$i\in\Ob(I)$.
\end{enumerate}
\end{enumerate}
\end{proposition}
\begin{proof}(i) follows immediately from corollary
\ref{cor_mitchell}(i,ii) and proposition \ref{prop_soberify}.

(ii.a): By virtue of proposition
\ref{prop_filter-Deligne}, we may assume that
$I$ is a partially ordered set, and to ease
notation, for every $u:i\to j$ in $I$ ({\em i.e.}
for every $i,j\in I$ such that $i\leq j$), we set
$p_{ij}:=p_u$.

Now, consider the family $\cF$ consisting of all
compatible systems $Z_\bullet:=(Z_i~|~i\in I)$, where :
\begin{itemize}
\item
$Z_i$ is a non-empty closed subset of $X_i$, for every $i\in I$
\item
$p_{ij}(Z_i)\subset Z_j$ for every $i,j\in I$ with $i\leq j$
\end{itemize}
and endow $\cF$ with a partial ordering, by declaring
that, for any $Z_\bullet,Z'_\bullet\in\cF$, we have
$Z_\bullet\leq Z'_\bullet$ if and only if $Z_i\subset Z'_i$
for every $i\in I$. Suppose now that
$(Z_\bullet^\lambda~|~\lambda\in\Lambda)$ is a totally
ordered subset of $\cF$, and set
$C_i:=\bigcap_{\lambda\in\Lambda}Z^\lambda_i$ for every $i\in I$.

\begin{claim}\label{cl_minimal-filter}
The resulting system $C_\bullet$ lies in $\cF$.
\end{claim}
\begin{pfclaim} Indeed, the assertion comes down to
saying that $C_i\neq\emptyset$ for every $i\in I$.
But notice that, for every finite subset
$\Lambda'\in\Lambda$ and every $i\in I$, the intersection
$\bigcap_{\lambda\in\Lambda'}Z_i^\lambda$ is obviously
closed and non-empty in $X_i$. Since $X_i$ is
quasi-compact for every $i\in I$, the claim follows.
\end{pfclaim}

Since $(X_i~|~i\in I)\in\cF$, we have $\cF\neq\emptyset$.
From claim \ref{cl_minimal-filter} and Zorn's lemma,
we deduce that $\cF$ admits minimal elements. We notice:

\begin{claim}\label{cl_dense-irred}
Let $Z_\bullet$ be a minimal element of $\cF$. We have :
\begin{enumerate}
\item
$Z_i$ is irreducible in $X_i$, for every $i\in I$.
\item
$p_{ij}(Z_i)$ is dense in $Z_j$, for every $i,j\in I$
with $i\leq j$.
\end{enumerate}
\end{claim}
\begin{pfclaim}(i): Suppose, by way of contradiction, that
there exist $i\in I$ and closed non-empty proper subsets
$Z^{(1)},Z^{(2)}$ of $Z_i$, such that $Z_i=Z^{(1)}\cup Z^{(2)}$.
For $t=1,2$, we consider the compatible system
$Z^{(t)}_\bullet$ such that
$$
Z^{(t)}_j:=\left\{\begin{array}{lll}
              Z_j\cap p^{-1}_{ji}Z^{(t)} & \qquad & \text{if $j\leq i$} \\
              Z_j & \qquad & \text{otherwise}.
                           \end{array}\right.
$$
By the minimality of $Z_\bullet$, neither $Z^{(1)}_\bullet$
nor $Z_\bullet^{(2)}$ lies in $\cF$, hence there must
exist $j,k\in I$ such that $Z^{(1)}_j=Z^{(2)}_k=\emptyset$.
Since $I$ is cofiltered, we may then assume that $j=k$,
in which case
$Z_k=Z_k\cap p^{-1}_{ki}Z_i=Z^{(1)}_k\cup Z^{(2)}_k=\emptyset$,
which is absurd.

(ii): Fix any $i,j\in I$ with $i\leq j$, and denote by
$T$ the topological closure of $p_{ij}(Z_i)$ in $Z_j$.
We consider the compatible system $T_\bullet$ such that
$$
T_k:=\left\{\begin{array}{lll}
              p^{-1}_{kj}T & \qquad & \text{if $k\leq j$} \\
              Z_k & \qquad & \text{otherwise}. \\
             \end{array}\right.
$$
Notice that $T_k=Z_k$ for every $k\leq i$; especially
$T_k\neq\emptyset$ for every $k\in I$, and therefore
$T_\bullet\in\cF$. By the minimality of $Z_\bullet$, we
must have $T_\bullet=Z_\bullet$, and therefore $T=Z_j$,
as stated.
\end{pfclaim}

Let now $Z_\bullet$ be any minimal element of $\cF$; according
to claim \ref{cl_dense-irred}, every $Z_i$ is irreducible,
say with generic point $\eta_i$, and $p_{ij}$ maps $\eta_i$
to $\eta_j$, whenever $i,j\in I$ with $i\leq j$. The
compatible system of points $(\eta_i~|~i\in I)$ is then 
a point of $X$, so the latter is not empty.

Next, notice that every open covering of $X$ can be
refined to a covering of the form :
\set\begin{equation}\label{eq_cofinal-cover}
X=\bigcup_{i\in I}p^{-1}_iU_i
\qquad
\text{for a family of open subsets $U_i\in X_i$}
\end{equation}
hence, in order to prove that $X$ is quasi-compact,
it suffices to check that, for a covering as
\eqref{eq_cofinal-cover}, there exists a finite
subset $J\subset I$ such that the family
$(p^{-1}_jU_j~|~j\in J)$ already covers $X$. To this
aim, set
$$
Z_i:=X_i\!\setminus\!U_i
\qquad\text{and}\qquad
Z'_i:=\bigcap_{j\geq i}p_{ij}^{-1}Z_j
\qquad
\text{for every $i\in I$}.
$$
Notice that \eqref{eq_cofinal-cover} means that
$\bigcap_{i\in I}p_i^{-1}Z_i=\emptyset$, and {\em a fortiory}
we have
\set\begin{equation}\label{eq_aswell}
\bigcap_{i\in I}p_i^{-1}Z'_i=\emptyset
\end{equation}
as well. However, by construction we have
$p_{ij}(Z'_i)\subset Z'_j$ for every $i,j\in I$ with
$i\leq j$, so $(Z'_i~|~i\in I)$ is a cofiltered system
of sober and quasi-compact topological spaces (lemma
\ref{lem_irreducible}(iii)), and \eqref{eq_aswell}
shows that the limit of this system is the empty
topological space; in light of the foregoing, we
conclude that there exists $i\in I$ such that
$Z'_i=\emptyset$. Since $X_i$ is quasi-compact,
this in turns means that there exists a finite
subset $J\subset\{j\in I~|~j\geq i\}$ such that
$\bigcap_{j\in J}p_{ij}^{-1}Z_j=\emptyset$. Lastly, this
means precisely that the family $(p^{-1}_jU_j~|~j\in J)$
covers $X$.

(ii.b): Fix $i_0\in I$, and let $U_{i_0}$ be any quasi-compact
open subset of $X_{i_0}$; the subset $I_0:=\{j\in I~|~j\leq i_0\}$
is cofinal in $I$, so we may replace $I$ by $I_0$ and
assume that $i_0$ is the largest element of $I$. In this
case, for every $j\leq i_0$ we may set $U_j:=p^{-1}_{ji_0}U$,
and $(U_i~|~i\in I)$ is a cofiltered system of quasi-compact
and sober topological spaces (lemma
\ref{lem_irreducible}(iii)), hence its limit $U$ is
quasi-compact (and sober), by (i) and (ii.a); on the other
hand, the induced map $U\to X$ is an open immersion, whence
the contention.
\end{proof}

\begin{remark}\label{rem_Hochster}
(i)\ \
Proposition \ref{prop_lim-qc-sober} fails in general
for non-cofiltered systems of quasi-compact and sober
spaces. For instance, a fibre product of such spaces
shall not be quasi-compact, in general.

(ii)\ \
It is clear that the topological space underlying every
affine scheme is spectral, and that of any scheme is
locally spectral. Conversely, it is shown in \cite{Hoch2}
that to every spectral space $X$ one can attach naturally
a commutative ring $A_X$ whose spectrum is isomorphic to
$X$. This rule $X\mapsto A_X$ is functorial with respect
to continuous maps of a certain type, which are introduced
in the next definition.
\end{remark}

\begin{definition}\label{def_constructible}
Let $f:X\to Y$ be a continuous map of topological
spaces, $T\subset X$ any subset, and endow $T$ with
the topology induced from $X$.
\begin{enumerate}
\item
We say that $f$ is {\em quasi-compact} (resp.
{\em quasi-separated}) if, for every quasi-compact
(resp. quasi-separated) open subset $U\subset X$, the
preimage $f^{-1}U$ is quasi-compact (resp. quasi-separated).
\item
We say that $f$ is {\em spectral} if, for every
pair of quasi-compact quasi-separated open subsets
$U\subset X$, $V\subset Y$ with $f(U)\subset V$,
the restriction $f_{|U}:U\to V$ is quasi-compact.
\item
We say that $T$ is {\em retro-compact in $X$} if the
inclusion map $T\to X$ is quasi-compact.
\item
We say that $T$ is {\em globally constructible} if
$T$ lies in the boolean algebra of subsets of $X$
generated by all retro-compact open subsets of $X$.
\item
We say that $T$ is {\em constructible} if every point
of $X$ admits an open neighborhood $U$ such that
$T\cap U$ is globally constructible in $U$.
\item
We say that $T$ is {\em pro-constructible} (resp.
{\em ind-constructible}) if every point of $X$
admits an open neighborhood $U$ such that $T\cap U$
is the intersection (resp. the union) of constructible
subsets of $U$.
\end{enumerate}
\end{definition}

\begin{remark}\label{rem_strongly-spectral}
Let $f:X\to Y$ and $g:Y\to Z$ be two continuous maps
of topological spaces.

(i)\ \
If $f$ and $g$ are quasi-compact (resp. quasi-separated),
then the same holds for $g\circ f$.

(ii)\ \
If $f$ and $g$ are spectral, and both $X$ and $Y$ are
locally coherent, then $g\circ f$ is spectral. Indeed,
say that $U\subset X$ and $V\subset Z$ are quasi-compact
and quasi-separated open subsets such that
$g\circ f(U)\subset V$; for every $x\in U$ we may find
a quasi-compact and quasi-separated open subset
$W_x\subset Y$ such that $g(W_x)\subset V$ and
$f(x)\in W_x$. Since $X$ is locally coherent, we may
find a quasi-compact and quasi-separated open neighborhood
$U_x$ of $x$ in $X$, such that $f(U_x)\subset W_x$.
Then both $f_{|U_x}:U_x\to W_x$ and $g_{|W_x}:W_x\to V$
are quasi-compact, so (i) says that the same holds for
their composition $(g\circ f)_{|U_x}:U_x\to V$. However,
finitely many such $U_x$ suffice to cover $U$, whence
the assertion.

(iii)\ \
Let us say that $f$ is {\em strongly spectral} if the
following holds. For every quasi-compact inclusion map
$U\to U'$ between open subsets of $Y$, the induced
open immersion $f^{-1}U\to f^{-1}U'$ is also quasi-compact.
We claim that if $f$ is strongly spectral, then $f$ is
spectral. Indeed, let $U\subset X$, $V\subset Y$ be two
quasi-compact and quasi-separated open subsets, and
$V'\subset V$ another quasi-compact open subset; by
assumption, the induced inclusion map $i:f^{-1}V'\to f^{-1}V$
is quasi-compact, so $f^{-1}_{|U}V'=i^{-1}U$ is quasi-compact,
whence the claim.

(iv)\ \
Obviously, a composition of strongly spectral maps
is strongly spectral. Moreover, it is easily seen that
any open immersion $i:U\to X$ is strongly spectral.
In view of (iii), it follows that every open immersion
is spectral.
\end{remark}

\begin{lemma}\label{lem_charact-coh-sob-sp}
Let $X$ be any topological space. We have :
\begin{enumerate}
\item
$X$ is quasi-separated if and only if it admits a covering
consisting of retro-compact open quasi-separated subsets.
\item
$X$ is sober if and only if it admits a covering consisting
of sober open subsets.
\item
The following conditions are equivalent :
\begin{enumerate}
\item
$X$ is locally coherent and quasi-separated.
\item
$X$ admits a basis consisting of quasi-compact open subsets,
which is closed under finite intersections.
\end{enumerate}
\item
The following conditions are equivalent :
\begin{enumerate}
\item
$X$ is coherent (resp. spectral).
\item
$X$ is quasi-compact and admits a covering consisting of
coherent (resp. spectral) retro-compact open subsets.
\item
$X$ is quasi-compact, quasi-separated and locally coherent
(resp. locally spectral).
\end{enumerate}
\item
$X$ is noetherian if and only if it is locally noetherian
and quasi-compact.
\item
If $X$ is locally noetherian, then $X$ is locally coherent
and quasi-separated.
\end{enumerate}
\end{lemma}
\begin{proof}(i): The condition is trivially necessary.
Suppose then that $(U_i~|~i\in I)$ is a covering of $X$
such that $U_i$ is quasi-separated, retro-compact and
open in $X$ for every $i\in I$; pick any two quasi-compact
open subsets $U$, $U'$. We may find a finite subset
$J\subset I$ such that $U\cup U'\subset\bigcup_{j\in J}U_j$,
and our assumptions imply that, for every $j\in J$, the
subset $V_j:=U\cap U'\cap U_j$ is quasi-compact, so the
same holds for $\bigcup_{j\in J}V_j=U\cap U'$.

(ii): Again, necessity is trivial. For the converse,
let $(U_i~|~i\in I)$ be an open covering of $X$ such
that $U_i$ is sober for every $i\in I$. Say that
$Z\subset X$ is any irreducible closed subset, and
that $Z_i:=U_i\cap Z\neq\emptyset$ for some $i\in I$;
then $Z_i$ is an irreducible subset of $U_i$ and
the topological closure of $Z_i$ in $X$ equals $Z$
(lemma \ref{lem_irreducible}(ii)). Let $\eta$
be the generic point of $Z_i$; then clearly the
topological closure of $\{\eta\}$ in $X$ equals $Z$.
It remains to check that $\eta$ is the unique element
of $Z$ with this property. However, say that $\eta'$
is another such element; then it is easily seen that
$\eta'\in U_i$ as well, and $\eta'$ is then also the
generic point of $Z_i$, whence the contention.

(iii): It is easily seen that (b)$\Rightarrow$(a).
Conversely, suppose that (a) holds, and choose a
covering $(U_i~|~i\in I)$ of $X$ consisting of
coherent open subsets; for each $i\in I$ we then
have a basis $(U_{ij}~|~j\in J_i)$ of the topology
of $U_i$, consisting of quasi-compact open subsets.
Set $\cF:=(U_{ij}~|~i\in I,\ j\in I_j)$, and let $\cF'$
be the family of all finite non-empty intersections
of elements of $\cF$; then $\cF'$ is a basis of
the topology of $X$, and since $X$ is quasi-separated,
any element of $\cF'$ is still a quasi-compact open
subset of $X$, whence (b). 

(iv): Obviously (a)$\Rightarrow$(b), and
(b)$\Rightarrow$(c) by virtue of (i) and remark
\ref{rem_sorite-spectral}(i). Lastly,
(c)$\Rightarrow$(a), by (iii) and (ii).

(v): Suppose that $X$ is locally noetherian and
quasi-compact, in which case it admits a finite
covering $(U_i~|~i\in I)$ consisting of open noetherian
subsets. Now, let $Z_\bullet:=(Z_n~|~n\in\N)$ be
a descending chain of closed subsets of $X$; by
assumption, every resulting chain
$(Z_n\cap U_i~|~n\in\N)$ is stationary, and since
$I$ is finite, it follows that $Z_\bullet$ is stationary
as well.

(vi): $X$ is locally coherent, by virtue of remark
\ref{rem_sorite-spectral}(iv); next, if $U,V\subset X$
are quasi-compact open subsets, then they are noetherian,
by (v), so their intersection is quasi-compact.
\end{proof}

\begin{proposition}\label{prop_local-qc-qs}
Let $f:X\to Y$ be a continuous map of topological
spaces, $(Y_i~|~i\in I)$ an open covering of\/ $Y$,
and $f_i:f^{-1}Y_i\to Y_i$ the restriction of $f$,
for every $i\in I$. We have :
\begin{enumerate}
\item
If\/ $Y$ is locally coherent, then $f$ is quasi-compact
if and only if $f_i$ is quasi-compact for every $i\in I$.
\item
If both $X$ and $Y$ are locally coherent, then $f$ is
quasi-separated and spectral if and only if $f_i$ is
quasi-separated and spectral for every $i\in I$. 
\item
Suppose that both $X$ and $Y$ are locally coherent, and
let $(X_i~|~i\in I)$ be an open covering of $X$ such that
$X_i\subset f^{-1}Y_i$ for every $i\in I$. Then $f$ is
spectral if and only if the restriction $X_i\to Y_i$
of $f$ is spectral for every $i\in I$.
\end{enumerate}
\end{proposition}
\begin{proof}(i): Clearly, if $f$ is quasi-compact, the
same holds for every $f_i$. Hence we may assume that
$f_i$ is quasi-compact for every $i\in I$, and we have
to check that $f$ is quasi-compact. To this aim, let
$U\subset Y$ be any quasi-compact open subset, and
set $U_i:=Y_i\cap U$ for every $i\in I$; since $Y$
is locally coherent, we may find, for every $i\in I$,
a covering $(U_{i,\lambda}~|~\lambda\in\Lambda_i)$ of
$U_i$ consisting of quasi-compact open subsets.
Since $U$ is quasi-compact, finitely many of these
$U_{i,\lambda}$ suffice to cover $U$, hence $f^{-1}U$
is covered by finitely many subsets of the type
$f^{-1}_iU_{i,\lambda}$, and the latter are quasi-compact
by assumption, so the same holds for $U$.

(iii): To begin with, we remark

\begin{claim}\label{cl_strongly-spectral}
Suppose that both $X$ and $Y$ are locally coherent.
Then $f$ is spectral if and only if it is strongly
spectral (see remark \ref{rem_strongly-spectral}(iii)).
\end{claim}
\begin{pfclaim} In view of remark \ref{rem_strongly-spectral}(iii),
we may assume that $f$ is spectral, and we show that $f$
is strongly spectral. Indeed, let $U\subset U'$ be any
inclusion of open subset of $Y$, with $U$ retro-compact
in $U'$. We need to check that $f^{-1}U$ is retro-compact
in $f^{-1}U'$, and by (i), the assertion can be checked
locally on $f^{-1}U'$. However, since $X$ and $Y$ are
locally coherent, for every $x\in f^{-1}U'$ we may find
quasi-compact and quasi-separated open neighborhoods $V'$
of $x$ in $f^{-1}U'$ and $W'$ of $f(x)$ in $U'$ such that
$f(V')\subset W'$. By assumption, $W'\cap U$ is quasi-compact,
so the same holds for $V:=V'\cap f^{-1}U$ (because $f$ is
spectral); since $V'$ is quasi-separated, the inclusion map
$V\to V'$ is then quasi-compact, as required.
\end{pfclaim}

Now, in light of claim \ref{cl_strongly-spectral}, for
the proof of (iii) we may assume that each restriction
$X_i\to Y_i$ is strongly spectral, and it suffices to
check that the same holds for $f$. However, let $U_1\subset U_2$
be any quasi-compact inclusion map of open subsets of $Y$;
by (i), the induced inclusion map $U_1\cap Y_i\to U_2\cap Y_i$
is also quasi-compact, so by assumption the same holds
for the inclusion
$$
X_i\cap f^{-1}U_1=X_i\cap f^{-1}(U_1\cap Y_i)\to
X_i\cap f^{-1}(U_2\cap Y_i)=X_i\cap f^{-1}U_2
\qquad
\text{for every $i\in I$}.
$$
Again (i) then implies that the inclusion $f^{-1}U_1\to f^{-1}U_2$
is quasi-compact, as required.

(ii): Again, we may assume that $f_i$ is quasi-separated
and spectral for every $i\in I$, and we show that the same
follows for $f$. By (iii), we know already that $f$ is
spectral, so we let $U\subset Y$ be any quasi-separated open
subset, and show that $f^{-1}U$ is quasi-separated in $X$.
To this aim, choose a family of quasi-compact open subsets
$(U_{i,\lambda}~|~\lambda\in\Lambda_i)$ which gives a covering
of $U\cap Y_i$, for every $i\in I$; since $U$ is quasi-separated,
$U_{i,\lambda}$ is retro-compact in $U$ for every $i\in I$ and
every $\lambda\in\Lambda_i$, hence $f^{-1}U_{i,\lambda}$ is
retro-compact in $f^{-1}U$, by claim \ref{cl_strongly-spectral}.
Now, let $V_1,V_2\subset f^{-1}U$ be any two quasi-compact
open subsets; it follows that $V_1\cap f^{-1}U_{i,\lambda}$
and $V_2\cap f^{-1}U_{i,\lambda}$ are quasi-compact for every
$i\in I$ and every $\lambda\in\Lambda_i$. Since $f_i$ is
quasi-separated, we deduce that
$V_{i,\lambda}:=V_1\cap V_2\cap f^{-1}U_{i,\lambda}$ is also
quasi-compact, for every such $i$ and $\lambda$.
Since $V_1$ is quasi-compact, finitely many of these
$f^{-1}U_{i,\lambda}$ suffice to cover $V_1$, and therefore
finitely many $V_{i,\lambda}$ cover $V_1\cap V_2$, so the
latter is quasi-compact, as required.
\end{proof}

\begin{example} We return briefly to remark
\ref{rem_Hochster}(ii) : with the terminology of
definition \ref{def_constructible}(ii), we may now
say that the continuous map $\Spec\,\phi:\Spec\,B\to\Spec\,A$
attached to any ring homomorphism $\phi:A\to B$ is spectral
(verification left to the reader). In light of
proposition \ref{prop_local-qc-qs}(iii) it follows that
the continuous map underlying any morphism of schemes
$f:X\to Y$ is spectral as well. Moreover, let $(Y_i~|~i\in I)$
be any covering of $Y$ consisting of open subschemes;
proposition \ref{prop_local-qc-qs}(i,ii) implies that
$f$ is quasi-compact (resp. quasi-separated) if and only
the same holds for each restriction $f_{|Y_i}:f^{-1}Y_i\to Y_i$.
\end{example}

\begin{lemma}\label{lem_sorite-construct}
Let $X$ be a topological space, $U\subset X$
an open subset, $T\subset X$ any subset, and endow
$U$ and $T$ with the topology induced from $X$.
The following holds :

{\em(i)}\ \
Every closed subset of $X$ is retro-compact. The union
of any finite number of retro-compact subsets of $X$ is
still retro-compact. The intersection of any open
retro-compact subset of $X$ with an arbitrary retro-compact
subset of $X$ is retro-compact.

{\em(ii)}\ \
Every globally constructible subset of $X$ is
retro-compact, and it is a finite union of subsets
of the form $V\!\setminus\!V'$, where $V$ and $V'$
are retro-compact open subsets of $X$.

{\em(iii)}\ \
Suppose that $T$ is retro-compact, and let $S\subset T$
be any subset. If $S$ is retro-compact in $T$, then
it is retro-compact in $X$ as well. If\/ $T$ is also
open and $S$ is globally constructible in $T$, then $S$
is globally constructible in $X$ as well.

{\em(iv)}\ \
$T$ is pro-constructible if and only if
$X\!\setminus\!T$ is ind-constructible.

{\em(v)}\ \
If\/ $T$ is retro-compact (resp. globally constructible,
resp. pro-constructible, resp. ind-constructible) in
$X$, then $T\cap U$ is retro-compact (resp. globally
constructible, resp. constructible, resp.
pro-constructible, resp. ind-constructible) in $U$.

{\em(vi)}\ \
The constructible subsets of $X$ form a boolean algebra.

{\em(vii)}\ \
If\/ $X$ is quasi-separated and $U$ is quasi-compact,
then $U$ is retro-compact in $X$.

{\em(viii)}\ \
Suppose that $X$ is quasi-compact; then we have :
\begin{enumerate}
\alphaenu
\item
If\/ $T$ is retro-compact in $X$, then $T$ is quasi-compact.
\item
Every globally constructible subset of $X$ is quasi-compact.
\end{enumerate}

{\em(ix)}\ \
Suppose that $X$ is quasi-compact and quasi-separated;
then we have :
\begin{enumerate}
\alphaenu
\item
$U$ is retro-compact in $X$ if and only if it is
quasi-compact.
\item
The globally constructible subsets of $X$ are precisely
the finite unions of subsets of the form $V\!\setminus\!V'$,
where $V$ and $V'$ are arbitrary quasi-compact open subsets
of $X$.
\end{enumerate}

{\em(x)}\ \
Suppose that $X$ is coherent; then we have :
\begin{enumerate}
\alphaenu
\item
$T$ is constructible in $X$ if and only if it is
globally constructible in $X$ (in which case $T$ is
also retro-compact, by {\em(ii)}).
\item
$T$ is pro-constructible (resp. ind-constructible) if
and only if it is the intersection (resp. the union)
of a family of constructible subsets of $X$.
\item
The constructible open (resp. closed) subsets of $X$
are precisely the quasi-compact open subsets (resp.
the complements of quasi-compact open subsets) of $X$.
\item
If\/ $T$ is retro-compact in $X$ and $S$ is any constructible
(resp. pro-constructible, resp. ind-constructible) subset of
$X$, then $T\cap S$ is constructible (resp. pro-constructible,
resp. ind-constructible) in $T$.
\item
If\/ $T$ is retro-compact in $X$, then it is a coherent
topological space, and the inclusion map $i_T:T\to X$ is
quasi-separated.
\item
If\/ $T$ is constructible in $X$, and $S$ is any
constructible subset of\/ $T$, then $S$ is constructible
in $X$.
\end{enumerate}

{\em(xi)}\ \
If\/ $X$ is spectral and $T$ is constructible in $X$,
then $T$ is spectral as well.
\end{lemma}
\begin{proof}(i) shall be left to the reader, and
(ii) follows easily from (i).

(iii): The first assertion is obvious, and the second
follows from (ii) and the first assertion.

(v): The proof for the case where $T$ is retro-compact
or globally constructible is left to the reader.
Suppose that $T$ is constructible in $X$; then,
for any given point $x\in U$, we may find an open
neighborhood $U'$ of $x$ in $X$, such that $T\cap U'$
is globally constructible in $U'$; then by the foregoing
case, we know that $U'\cap U$ is an open neighborhood of
$x$ in $U$ such that $T\cap U'\cap U$ is globally
constructible in $U'\cap U$, so $T\cap U$ is constructible
in $U$. The assertion for the case where $T$ is
pro-constructible or ind-constructible, is an immediate
consequence.

(vi): It is easily seen that the complement of
a constructible subset is constructible. Next,
suppose that $T$ and $T'$ are two constructible
subsets of $X$, and $x\in X$ is any point; then
we may find open neighborhoods $U$ and $U'$ of
$x$ such that $T\cap U$ (resp. $T'\cap U'$) is
globally constructible in $U$ (resp. in $U'$).
By (v), it follows that $T\cap U\cap U'$ and
$T\cap U\cap U'$ are globally constructible in
$U\cap U'$, and therefore the same holds for both
$(T\cap T')\cap U\cap U'$ and $(T\cup T')\cap U\cap U'$.
We conclude that both $T\cap T'$ and $T\cup T'$
are constructible, whence the assertion.

(iv) follows easily from from (vi).

(vii) and (viii.a) are obvious, and (viii.b) follows
from (viii.a) and (ii).

(ix.a) follows from (vii) and (viii.a), and (ix.b)
follows from (ix.a) and (ii).

(x.a): Indeed, by assumption we may find for every
$x\in X$ an open neighborhood $U_x$ of $x$ in $X$
such that $T\cap U_x$ is globally constructible in
$U_x$; by (v) we may then assume that $U_x$ is
quasi-compact, so $U_x$ is retro-compact, by (ix.a)
and remark \ref{rem_sorite-spectral}(i), and
consequently $T\cap U_x$ is globally constructible
in $X$, by (iii). Finitely many of such $U_x$ suffice
to cover $X$, whence the assertion. 

(x.b): Indeed, say that $T$ is ind-constructible;
arguing as in the proof of (x.a) we may find, for
every $x\in X$, a quasi-compact open subset $U_x$
such that $T\cap U_x$ is the union of a family
$(T_{i,x}~|~i\in I_x)$ of constructible subsets of
$U_x$. By (x.a) and
remark \ref{rem_sorite-spectral}(ii), each $T_{i,x}$ is
globally constructible in $U_x$, hence also in $X$, by
(viii.a), (iii) and remark \ref{rem_sorite-spectral}(i);
since $T=\bigcup_{x\in X}\bigcup_{i\in I_x}T_{i,x}$,
the claim follows. The assertion for pro-constructible
subsets is reduced to this case, by taking complements
and using (iv).

(x.c): Indeed, every quasi-compact open subset is
constructible, by (vii); conversely, (x.a) and
(ix.b) show that every constructible subset of $X$
is quasi-compact. The assertion about constructible
closed subsets of $X$ then follows immediately,
taking into account (vi).

(x.d): In view of (x.b), it suffices to consider
the case where $S$ is constructible, and so $S$
is even globally constructible in $X$, by (x.a).
Then, we are further reduced to the case where $S$
is an open retro-compact subset of $X$, {\em i.e.}
$S$ is open and quasi-compact in $X$, by (x.c), and
it suffices to check that $T\cap S$ is retro-compact
in $T$. Thus, let $V$ be any quasi-compact open subset
of $T$; since $X$ is coherent, we may write $V=V'\cap T$,
for some quasi-compact open subset $V'$ of $X$ (details
left to the reader), and then
$V\cap(T\cap S)=(V'\cap S)\cap T$. Now, $V'\cap S$
is open and quasi-compact in $X$ (remark
\ref{rem_sorite-spectral}(i)) and therefore $V\cap(T\cap S)$
is quasi-compact in $T$, and therefore also in $T\cap S$,
since the latter is open in $T$.

(x.e): $T$ is quasi-compact by (viii.a), and if
$(U_i~|~i\in I)$ is a basis of quasi-compact open subsets
of $X$ which is closed under finite intersections, then
$(T\cap U_i~|~i\in I)$ is a basis of quasi-compact open
subsets for $T$ which is also closed under finite intersections,
so $T$ is coherent.
Next, let $U\subset X$ be any open subset (then $U$ is
quasi-separated, by remark \ref{rem_sorite-spectral}(i)),
and $V_1,V_2$ any two quasi-compact open subsets of
$U\cap T$; then $V_1$ and $V_2$ are also quasi-compact
in $T$, so the same follows for $V_1\cap V_2$ (again,
by remark \ref{rem_sorite-spectral}(i)), which shows
that $i_T$ is quasi-separated.

(x.f): In view of (x.e), the topological space $T$ is
coherent, so $S$ is globally constructible in $T$ by
(x.a), and in light of (ii) and (vi) we may assume that
$S=V\setminus V'$, for two retro-compact open subsets
$V,V'$ of $T$. Then $V$ and $V'$ are quasi-compact in
$T$, by (viii.a), and since $X$ is coherent, we may
find quasi-compact open subsets $W,W'$ of $X$ such
that $V=W\cap T$ and $V'=W'\cap T$, so that
$S=(W\setminus W')\cap T$, and the latter is
constructible in $X$, again by (vi) and (x.c).

(xi): In light of (x.e) and (x.a), it remains to show
that $T$ is sober, and taking into account (ii) and
lemma \ref{lem_irreducible}(iv), we may assume that
$T$ is locally closed in $X$, in which case it
suffices to invoke lemma \ref{lem_irreducible}(iii).
\end{proof}

\begin{remark}\label{rem_sorite-qcoh-maps}
Let $f:X\to Y$ be any quasi-compact continuous map
of topological spaces.

(i)\ \
If $T\subset X$ is a retro-compact subset of $X$,
then $f(T)$ is a retro-compact subset of $Y$. Indeed,
if $U\subset Y$ is any quasi-compact open subset, we
have $f(T)\cap U=f(T\cap f^{-1}U)$, and $f^{-1}U$ is
a quasi-compact open subset of $X$, by assumption.

(ii)\ \
If $f$ is quasi-separated, $Y$ is locally coherent,
and $T\subset Y$ is any constructible (resp.
pro-constructible, resp. ind-constructible) subset
in $Y$, then $f^{-1}T$ enjoys the same property in
$X$. To see this, we may assume that $Y$ is coherent,
in which case $X$ is quasi-separated and $T$ is
constructible (resp. an intersection, resp. a union
of constructible subsets, by lemma
\ref{lem_sorite-construct}(x.b)); then it suffices
to consider the case where $T$ is constructible,
and we are further reduced to the case where
$T$ is open and quasi-compact (lemma
\ref{lem_sorite-construct}(ix.b,x.a)), for which the
assertion follows from lemma \ref{lem_sorite-construct}(vii).

(iii)\ \
If $f$ is quasi-separated, then $f$ is spectral. Indeed,
say that $U\subset X$ and $V\subset U$ are quasi-compact
and quasi-separated open subsets such that $f(U)\subset V$,
and let $V'\subset V$ be any quasi-compact open subset of
$V$; then $f^{-1}V$ is quasi-compact and quasi-separated,
and $f_{|U}^{-1}V'=U\cap f^{-1}V'$, which is quasi-compact,
under the current assumptions, because $U$ and $f^{-1}V'$
are quasi-compact open subsets of $f^{-1}V$.

(iv)\ \
Suppose that $X$ is locally coherent, $T\subset X$ is any
retro-compact subset, and endow $T$ with the topology
induced from $X$. Then the inclusion map $i_T:T\to X$ is
spectral. To see this, first we invoke proposition
\ref{prop_local-qc-qs}(i,iii) to reduce to the case
where $X$ is coherent. Then the assertion follows from
(iii) and lemma \ref{lem_sorite-construct}(x.e). Together
with lemma \ref{lem_sorite-construct}(i), we see especially
that any closed immersion $Z\to X$ is quasi-compact,
quasi-separated and spectral.
\end{remark}

\begin{definition}\label{subsec_constr-top}
Let $(X,\cT)$ be a topological space. The set of
ind-constructible subsets of $X$ is closed under
finite unions and finite intersections, and forms
a basis for a topology $\cC$ on $X$ called the
{\em constructible topology}. We shall denote by
$X^c$ the topological space $(X,\cC)$.
\end{definition}

\begin{remark}\label{rem_constr-top}
Let $f:X\to Y$ be a continuous map of topological
spaces.

(i)\ \
The map $f$ is not necessarily continuous with respect
to the constructible topologies on $X$ and $Y$, but
remark \ref{rem_sorite-qcoh-maps}(ii) says that if
$f$ is quasi-separated and quasi-compact, and $Y$
is locally coherent, then the induced map
$$
f^c:X^c\to Y^c
$$
(whose underlying map of sets is the same as $f$)
is continuous.

(ii)\ \
Let $U$ be any open subset of a topological space $X$,
and endow $U$ with the topology induced from $X$. Then
the topology of $U^c$ is finer than the topology on $U$
induced from the topology of $X^c$ (lemma
\ref{lem_sorite-construct}(v)). The same holds if $X$
is coherent and $U$ is retro-compact (lemma
\ref{lem_sorite-construct}(x.d)).

(iii)\ \
Suppose that $X$ is locally coherent. By lemma
\ref{lem_sorite-construct}(iv, x.b), it follows that
the open (resp. closed) subsets of $X^c$ are the 
ind-constructible (resp. pro-constructible) subsets
of $X$.
\end{remark}

\begin{proposition}\label{prop_constr-top}
Let $f:X\to Y$ be any spectral map between locally coherent
topological spaces, $T\subset Y$ any subset. The following holds :
\begin{enumerate}
\item
$f^c$ is continuous.
\item
If\/ $T$ is constructible (resp. ind-constructible, resp.
pro-constructible) in $Y$, then $f^{-1}T$ is constructible
(resp. ind-constructible, resp. pro-constructible) in $X$.
\end{enumerate}
\end{proposition}
\begin{proof} Clearly, it suffices to check (ii).
However, let $T\subset Y$ be any subset as in (ii), and
pick a covering $(U_i~|~i\in I)$ of $Y$ consisting
of open coherent subsets; it suffices to check
that $f^{-1}(T\cap U_i)$ is constructible (resp.
ind-constructible, resp. pro-constructible) in
$f^{-1}U_i$ for every $i\in I$, so we may replace
$Y$ by any $U_i$, and assume that $Y$ is coherent,
in which case, lemma \ref{lem_sorite-construct}(x.b)
implies that it suffices to consider the case where
$T$ is constructible in $Y$. Then, $T$ shall be a
finite union of subsets of the form $T'\setminus T''$,
where $T''\subset T'$ are any two quasi-compact open
subsets of $Y$ (lemma \ref{lem_sorite-construct}(ix.b,x.a)),
and we may further reduce to the case where $T$ is
an open quasi-compact subset of $Y$.
Now, let $U$ be any coherent open subset of $X$;
it suffices to check that $W:=U\cap f^{-1}V$ is
quasi-compact in $U$ (lemma
\ref{lem_sorite-construct}(x.c)), and the latter
holds, since $f$ is spectral.
\end{proof}

\begin{lemma}\label{lem_main-spectral}
Let $X$ be a locally coherent topological space,
$T\subset X$ any subset, and endow $T$ with the
topology induced from $X$. We have :
\begin{enumerate}
\item
If\/ $T$ is either constructible or open in $X$, the
topology on $T$ induced from $X^c$ agrees with
the topology of\/ $T^c$.
\item
The topology of $X^c$ is finer than that of $X$.
\end{enumerate}
\end{lemma}
\begin{proof}(i): It suffices to show that :
\begin{enumerate}
\alphaenu
\item
if $S$ is any ind-constructible subset of $X$, then
$S\cap T$ is ind-constructible in $T$
\item
if $S$ is any ind-constructible subset of $T$, then
$S$ is ind-constructible in $X$.
\end{enumerate}
However, fix a basis $(U_i~|~i\in I)$ of $X$
consisting of coherent open subsets, and set
$T_i:=T\cap U_i$ for every $i\in I$. If $S$ is
as in (a), and if $S\cap T_i$ is ind-constructible in
$T_i$ for every $i\in I$, then clearly $S\cap T$
is ind-constructible in $T$. If $S$ is as in (b), and
$S\cap U_i$ is ind-constructible in $U_i$ for every
$i\in I$, then $S$ is ind-constructible in $X$. Hence,
it suffices to prove (a) and (b) with $X$ and $T$
replaced by respectively $U_i$ and $T_i$, for every
$i\in I$. We may then assume from start that $X$ is
coherent. In this case, $T$ is either open or
retro-compact in $X$ (lemma \ref{lem_sorite-construct}(x.a)),
hence the topology of $T^c$ is finer than the topology
$\cT$ on $T$ induced from $X^c$
(remark \ref{rem_constr-top}(ii)), whence (a).
To show (b), suppose first that $T$ is constructible
in $X$; then, any ind-constructible subset of $T$ is
a union of constructible subsets in $T$ (lemma
\ref{lem_sorite-construct}(x.b,e)), and therefore it
is also a union of constructible subsets in $X$ (lemma
\ref{lem_sorite-construct}(x.f)), whence (b),
in this case. Lastly, if $T$ is open in $X$, and
$S$ is ind-constructible in $T$, pick a covering
$(T_i~|~i\in I)$ of $T$ consisting of open coherent
subsets; then $S_i:=S\cap T_i$ is ind-constructible
in $T_i$ for every $i\in I$, so for every such $i$
we may find a family $(S_{i,\lambda}~|~\lambda\in\Lambda_i)$
of globally constructible subsets of $T_i$, whose union
is $S_i$ (lemma \ref{lem_sorite-construct}(x.a,b)).
However, $T_i$ is retro-compact in $X$, so $S_{i,\lambda}$
is constructible in $X$, for every $i\in I$ and
every $\lambda\in\Lambda_i$ (lemma
\ref{lem_sorite-construct}(iii)); finally,
$S=\bigcup_{i\in I}\bigcup_{\lambda\in\Lambda_i}S_{i,\lambda}$
is ind-constructible in $X$, so (b) holds also in this
case.

(ii): We have to check that every open subset $U$
of $X$ is ind-constructible in $X$; however, $U$
is ind-constructible in $U$, so the assertion follows
from the foregoing.
\end{proof}

\sset\subsubsection{}\label{subsec_filter}
We come now to the main theorem of this section.
To prepare the proof, we introduce the following
terminology. Let $\cF$ be any family of subsets
of a topological space $T$; we say that $\cF$ is
a {\em closed filter} (resp. a
{\em constructible filter}) if :
\begin{itemize}
\item
every $Z\in\cF$ is closed (resp. constructible) in $T$
\item
$\emptyset\notin\cF$ and $\cF\neq\emptyset$
\item
if $Z,Z'\in\cF$, then $Z\cap Z'\in\cF$
\item
if $Z\in\cF$ and $Z'$ is any closed (resp. constructible)
subset of $T$ containing $Z$, then $Z'\in\cF$.
\end{itemize}
The {\em center} of a filter $\cF$ is the subset
$\bigcap_{Z\in\cF}Z$. It is easily seen that a topological
space $T$ is quasi-compact if and only if every closed
filter of $T$ has non-empty center.

\begin{theorem}\label{th_main-spectral}
Let $X$ be a topological space. We have :
\begin{enumerate}
\item
If $X$ is spectral, $X^c$ is compact.
\item
If $X$ is locally spectral and quasi-separated, then
$X^c$ is locally compact and separated.
\end{enumerate}
\end{theorem}
\begin{proof} (We say that a topological space is
{\em locally compact} if it admits a covering consisting
of compact open subsets.)

(i): We check first that $X^c$ is separated.
To this aim, let $x,y$ be any two points of $X$,
and denote by $\bar{\{x\}}$ and $\bar{\{y\}}$
the topological closures of $\{x\}$ and $\{y\}$
with respect to $\cT$. If $x\notin\bar{\{y\}}$, it
follows from lemma \ref{lem_main-spectral}(ii) that
we may find a constructible subset $T\subset X$ such
that $x\in T$ and $y\notin T$; then the complement
$S$ of\/ $T$ in $X$ is also constructible, and $T$
and $S$ separate the points $x$ and $y$. Likewise
one argues in the symmetric case where
$y\notin\bar{\{x\}}$. Lastly, if neither of these
conditions hold, we must have $x=y$, since $X$ is
sober and both $x$ and $y$ are generic points of
the same closed subset of $X$.

To show that $X^c$ is quasi-compact, we check that
its closed filters have non-empty center. However,
let $\cF$ be any such filter, and denote by $\cF'$
the set of all constructible subsets $T$ of $X$
such that there exists $F\in\cF$ with $F\subset T$.
Clearly $\cF'$ is a constructible filter of $X$,
thus, we come down to :

\begin{claim}(i)\ \
$\cF$ and $\cF'$ have the same center.

(ii)\ \
Every constructible filter of $X$ has non-empty center.
\end{claim}
\begin{pfclaim}(i): Clearly the center of $\cF'$ contains
the center of $\cF$. For the converse, it suffices to
remark that every closed subset of $X$ is pro-constructible
in $X$ (lemma \ref{lem_main-spectral}(ii) and remark
\ref{rem_constr-top}(iii)), and therefore it is the
intersection of the constructible subsets that contain
it (lemma \ref{lem_sorite-construct}(x.b)).

(ii): In view of lemma \ref{lem_sorite-construct}(xi),
every such filter can be regarded as a cofiltered system
of spectral spaces, with continuous transition maps.
Then the assertion is a special case of proposition
\ref{prop_lim-qc-sober}(ii.a).
\end{pfclaim}

(ii): Lemma \ref{lem_main-spectral} and (i) imply that
$X^c$ is locally compact. To see that $X$ is separated,
let $x,y\in X$ be any two distinct points. We may then
find a quasi-compact open subset $U$ of $X$ containing
both $x$ and $y$, and the induced map $U^c\to X^c$ is
an open immersion (lemma \ref{lem_main-spectral}), so
we may replace $X$ by $U$, and assume from start that
$X$ is is spectral (lemmata \ref{lem_charact-coh-sob-sp}(iv)
and \ref{lem_sorite-construct}(vii)), which is the case
covered by (i).
\end{proof}

\begin{corollary}\label{cor_hair-split}
Let $X$ be any spectral topological space, $U\subset X$
a subset. We have :
\begin{enumerate}
\item
$U$ is open and quasi-compact in $X^c$ if
and only if it is constructible in $X$.
\item
$X^c$ is a spectral topological space.
\item
If\/ $Y$ is any locally spectral topological space,
the same holds for $Y^c$.
\end{enumerate}
\end{corollary}
\begin{proof}(i): Suppose first that $U$ is open and
quasi-compact in $X^c$; then $U$ is a union of a
family $(U_i~|~i\in I)$ constructible subsets of
$X$ (lemma \ref{lem_sorite-construct}(x.b)), and
there exists a finite subset $J\subset I$ such that
$(U_i~|~i\in J)$ already covers $U$, so $U$ is
constructible.

Conversely, suppose that $U$ is constructible in $X$,
in which case it is both open and closed in $X^c$;
since the latter is quasi-compact (theorem
\ref{th_main-spectral}(i)), we see that $U$ is
quasi-compact in $X^c$.

(ii): Indeed, (i) and lemma \ref{lem_sorite-construct}(vi)
easily imply that $X^c$ is coherent, and $X^c$ is trivially
sober, since in any separated space the irreducible closed
subsets are exactly those that contain a unique point.

(iii) follows from (ii) and lemma
\ref{lem_main-spectral} : details left to the reader.
\end{proof}

\begin{corollary} Let $f:X\to Y$ be a continuous map
of locally spectral topological spaces. The following
holds :
\begin{enumerate}
\item
$f$ is spectral if and only if the same holds for $f^c$.
\item
Suppose furthermore that $f$ is spectral. Then :
\begin{enumerate}
\item
$f$ is quasi-compact if and only if the same holds
for $f^c:X^c\to Y^c$.
\item
$f$ is quasi-separated if and only if the same
holds for $f^c$.
\end{enumerate}
\end{enumerate}
\end{corollary}
\begin{proof} We notice first :

\begin{claim} In order to prove the corollary, we may
assume that $Y$ is spectral.
\end{claim}
\begin{pfclaim} Let $(Y_i~|~i\in I)$ be a covering
of $Y$ consisting of spectral open subsets; then
each $Y_i^c$ is open in $Y^c$, and its topology
agrees with the topology induced from $Y^c$ (lemma
\ref{lem_main-spectral}(i,ii)), so $(Y^c_i~|~i\in I)$
is an open covering of $Y^c$. Then the claim follows
immediately from proposition \ref{prop_local-qc-qs}(iii)
and corollary \ref{cor_hair-split}(iii).
\end{pfclaim}

Thus, we assume henceforth that $Y$ is spectral.
Next, we remark, quite generally :

\begin{claim}\label{cl_cont-is-spectral}
Let $g:T\to S$ be any continuous map between locally
spectral topological spaces, with $S$ locally compact.
Then $g$ is spectral.
\end{claim}
\begin{pfclaim} By lemma \ref{prop_local-qc-qs}(iii),
we may assume that $S$ is compact and that $T$ is
quasi-compact and quasi-separated.
Now, let $V\subset S$ be any quasi-compact open subset;
since $S$ is separated, $V$ is also closed in $S$, hence
$g^{-1}V$ is closed in $T$, so it is quasi-compact, since
the same holds for $T$. Thus, $g$ is quasi-compact, and
it is also obviously quasi-separated (due to remark
\ref{rem_sorite-spectral}(i)), so the claim follows
from remark \ref{rem_sorite-qcoh-maps}(iii).
\end{pfclaim}

(i): In light of claim \ref{cl_cont-is-spectral},
proposition \ref{prop_constr-top}(i) and theorem
\ref{th_main-spectral}(i), we see that if $f$ is
spectral, the same holds for $f^c$, hence we may
assume that $f^c$ is spectral, and we show that the
same holds for $f$. Thus, let $U\subset X$ and $V\subset Y$
be quasi-compact and quasi-separated open subsets
such that $f(U)\subset V$; it suffices to check that
$f_{|U}:U\to V$ is quasi-compact, and since $U^c$
(resp. $V^c$) is open in $X^c$ (resp. in $Y^c$, by
lemma \ref{lem_main-spectral}(ii)), the restriction
$f^c_{|U}:U^c\to V^c$ is still spectral, so we may
replace $X$ by $U$ and $Y$ by $V$, and assume that
$X$ is quasi-compact and quasi-separated, in which
case it is spectral (lemma \ref{lem_charact-coh-sob-sp}(iv)),
and we need to show that $f$ is quasi-compact.
Now, let $V\subset Y$ be any quasi-compact open subset;
then $V$ is spectral (remark \ref{rem_sorite-spectral}(ii)),
and $V^c$ is open and quasi-compact in $Y^c$ (lemma
\ref{lem_main-spectral}(ii) and theorem
\ref{th_main-spectral}(i)). Hence, $V^c$ is closed in
$Y^c$, and therefore $f^{-1}V^c$ is closed, and therefore
quasi-compact, in $X^c$. Since the topology of the latter
is finer than that of $X$, we conclude that $f^{-1}V$ is
quasi-compact in $X$.

(ii.a): If $f$ is quasi-compact, it follows that
$X$ admits a finite covering consisting of spectral
open subsets, and if $U$ is any such subset, the
restriction $f_{|U}$ is quasi-compact,
since $f$ is spectral. Clearly it suffices
to check that $f_{|U}^c$ is quasi-compact for
every such $U$, so we may replace $X$ by $U$ and
assume that $X$ is spectral as well, in which case,
notice that $f$ is obviously quasi-separated. Now,
say that $U\subset Y^c$ is open and quasi-compact;
so, $U$ is constructible in $Y$ (corollary
\ref{cor_hair-split}(i)), and then
$V:=f^{-1}U$ is constructible in $X$, by remark
\ref{rem_sorite-qcoh-maps}(ii), so it is
quasi-compact in $X^c$ (corollary \ref{cor_hair-split}(i)).

Conversely, say that $f^c$ is quasi-compact, and
$U$ is a quasi-compact open subset of $Y$. Then
$U$ is constructible in $Y$, so it is open and
quasi-compact in $Y^c$ (corollary \ref{cor_hair-split}(i)),
therefore $f^{-1}U$ is quasi-compact in $X^c$, and
since the topology of $X^c$ is finer than that of
$X$ (lemma \ref{lem_main-spectral}(ii)), we conclude
that $f^{-1}U$ is quasi-compact in $X$.

(ii.b): Suppose that $f$ is quasi-separated; then $X$
is quasi-separated, so $X^c$ is separated (theorem
\ref{th_main-spectral}(ii)), and it follows easily
that $f^c$ is quasi-separated (details left to the
reader). Conversely, if $f^c$ is quasi-separated,
then $X^c$ is quasi-separated, since $Y^c$ is
separated (by theorem \ref{th_main-spectral}(i)).
Then, say that $U,V$ are two quasi-compact open
subsets of $X$; we may cover $U$ (resp. $V$) by a
finite family $(U_i~|~i\in I)$ (resp. $(V_j~|~j\in J)$)
of spectral open subsets, and the resulting maps
$U_i^c\to X^c$, $V^c_j\to X^c$ are open immersions
(lemma \ref{lem_main-spectral}). Moreover, each
$U_i^c$ and $V^c_j$ is compact, by theorem
\ref{th_main-spectral}(i), so $U^c_i\cap V^c_j$ is
quasi-compact in $X^c$, for every $i\in I$ and every
$j\in J$; consequently $U\cap V$ is quasi-compact
in $X^c$, and {\em a fortiori} also in $X$, since
the topology of the latter is coarser than that
of $X^c$ (lemma \ref{lem_main-spectral}(ii)).
\end{proof}

\begin{corollary}\label{cor_main-spectral} 
{\em (i)}\ \
In the situation of \eqref{subsec_sober-lim}, suppose
moreover that :
\begin{enumerate}
\alphaenu
\item
The topological space $X_i$ is spectral for every $i\in\Ob(I)$.
\item
The transition map $p_u$ is quasi-compact, for every
morphism $u$ of $I$.
\end{enumerate}
\qquad\quad\ Then $X$ is spectral.
\begin{enumerate}
\addenu
\item
Suppose furthermore, that the system $X_\bullet$ is cofiltered.
Then :
\begin{enumerate}
\item
For every constructible (resp. open and quasi-compact)
subset $U\subset X$ there exist $i\in\Ob(I)$ and a constructible
(resp. open and quasi-compact) subset $U_i\subset X_i$
such that $U=p^{-1}_iU_i$.
\item
For any $i\in\Ob(I)$ and any two constructible subsets
$U,U'$ of $X_i$ such that $p^{-1}_iU\subset p^{-1}_iU'$,
there exists a morphism $u:j\to i$ in $I$, such
that $p^{-1}_uU\subset p^{-1}_uU'$.
\end{enumerate}
\end{enumerate}
\end{corollary}
\begin{proof}(i): By proposition
\ref{prop_lim-qc-sober}(i) we know already that $X$
is sober, so it remains only to check that $X$ is
coherent. However, under the stated conditions, the
datum $(X_i^c~|~i\in\Ob(I))$ is a system of quasi-compact
and separated topological spaces, with continuous
transition maps (theorem \ref{th_main-spectral} and
remark \ref{rem_constr-top}(i)); denote by $X^c$ the
limit of this system (in the category of topological
spaces). Now, the product $Y:=\prod_{i\in\Ob(I)}X^c_i$ is
likewise quasi-compact and separated, and we remark :

\begin{claim}\label{cl_alexandrov}
The natural map $X^c\to Y$ is a closed immersion.
\end{claim}
\begin{pfclaim} For every $i\in\Ob(I)$, denote by $q_i:Y\to X^c_i$
the natural projection; the image of $X^c$ in $Y$ is the
intersection of the subsets
$Y_u:=
\{y_\bullet\in Y~|~p_u\circ q_i(y_\bullet)=q_j(y_\bullet)\}$,
for every morphism $u:i\to j$ in $I$, so it suffices to
check that $Y_u$ is closed in $Y$, for every such $u$.
This is a special case of the following more general
assertion. Let $f,g:T\to T'$ be two continuous maps of
topological spaces, with $T'$ separated; then the subset
$$
S:=\{t\in T~|~f(t)=g(t)\}
$$
is closed in $T$. To show the latter, denote
$F:T\to T'\times T'$ the map such that $F(t)=(f(t),g(t))$
for every $t\in T$; then $S=F^{-1}\Delta_{T'}$, where
$\Delta_{T'}:=\{(t',t')~|~t'\in T'\}$ is the diagonal
of $T'\times T'$, and $F$ is continuous for the product
topology on $T'\times T'$. We then come down to showing
that $\Delta_{T'}$ is closed in $T'\times T'$, which
holds if (and only if) $T'$ is separated.
\end{pfclaim}

It follows from claim \ref{cl_alexandrov} that $X^c$
is quasi-compact and separated; the topology of $X^c$
is obviously finer than that of $X$, so the latter
is quasi-compact as well. It remains to check that $X$
admits a basis of quasi-compact open subsets closed
under finite intersections. To this aim, let
$J\subset\Ob(I)$ be any finite subset, pick a quasi-compact
open subset $U_j$ of $X_j$ for each $j\in J$, and set
$$
V:=\prod_{i\in I\setminus J}X_j\times\prod_{j\in J}U_j
\qquad
U:=X\cap V.
$$
It suffices to show that $U$ is quasi-compact in $X$.
However, notice that $U_j$ is constructible in $X_j$,
so it is both open and closed in $X^c_j$, and therefore
$V$ is closed in $Y$; since $X^c$ is closed in $Y$ by
claim \ref{cl_alexandrov}, it follows that $U$ is
quasi-compact in $X^c$, and {\em a fortiori}, also in $X$.

(ii.a): Suppose first that $U$ is open and quasi-compact;
then $U$ is a finite union of open subsets of the form
$U_J:=\bigcap_{j\in J}p_j^{-1}U_j$, where $J\subset\Ob(I)$ is
a finite subset, and $U_j$ is an open quasi-compact
subset of $X_j$, for every $j\in J$. For such a given
$U_J$, since $I$ is cofiltered, there exists $i\in\Ob(I)$
with morphisms $u_j:i\to j$ for every $j\in J$, and
therefore $U_J=p^{-1}_iV_i$, where
$V_i:=\bigcap_{j\in J}u_j^{-1}U_j$ is quasi-compact in
$X_i$, since each $u_j$ is quasi-compact and $X_i$ is
quasi-separated (remark \ref{rem_sorite-spectral}(ii)).
Thus, $U=\bigcup_{j\in J'}V_j$ for a finite subset
$J'\subset\Ob(I)$ and quasi-compact open subsets
$V_j\subset X_j$; again, we may find $i\in\Ob(I)$ and
a morphism $v_j:i\to j$ for every $j\in J'$, so that
$U=p^{-1}_iU_i$, where $U_i:=\bigcup_{j\in J'}v^{-1}_jV_j$
is open and quasi-compact in $X_i$.

The case where $U$ is constructible is reduced easily
to the foregoing case, taking into account (i) and
lemma \ref{lem_sorite-construct}(ix.b,x.a), and
arguing as above, using the assumption that  $X_\bullet$
is cofiltered (details left to the reader).

(ii.b)\ \
Let $Z:=U\!\setminus\!U'$, set $Z_u:=p^{-1}_uZ$
for every morphism $u:j\to i$ in $I$, and endow
$Z_u$ with the topology induced from $X_j^c$. From
theorem \ref{th_main-spectral}(i) and lemma
\ref{lem_sorite-construct}(xi) we see that
$Z_\bullet:=(Z_u~|~u\in\Ob(I/i))$ is a cofiltered
system of compact topological spaces (see example
\ref{ex_filtered-final}(i)), and the assumption
means that the limit of $Z_\bullet$ is empty, so there
exists $u\in\Ob(I/i)$ such that $Z_u=\emptyset$,
whence the assertion.
\end{proof}

\begin{corollary}\label{cor_procon-is-spec}
Let $X$ be any spectral topological space, $T\subset X$
a pro-constructible subset, and endow $T$ with the
topology induced from $X$. Then $T$ is spectral and
retro-compact in $X$. 
\end{corollary}
\begin{proof} By lemma \ref{lem_sorite-construct}(x.b),
we may write $T$ as the intersection of a system
$(T_i~|~i\in I)$ of constructible subsets of $X$.
Let $J$ be the set of all finite subsets of $I$,
and for every $S\in J$, endow $T_S:=\bigcap_{i\in S}T_i$
with the topology induced from $X$. Then $T_S$ is a
constructible subset of $X$ for every $S\in J$, and
clearly the topological space $T$ is the limit of the
(cofiltered) inverse system $(T_S~|~S\in J)$. By lemma
\ref{lem_sorite-construct}(xi) each such $T_S$ is spectral and
retro-compact in $X$ (lemma \ref{lem_sorite-construct}(x.a));
therefore, if $S\subset S'$ are any two elements of
$J$, then $T_S$ is constructible in $T_{S'}$ (lemma
\ref{lem_sorite-construct}(x.d)), hence also retro-compact
in $T_{S'}$ (lemma \ref{lem_sorite-construct}(x.a)),
{\em i.e.} the induced map $T_S\to T_{S'}$ is quasi-compact,
and the assertion follows from corollary
\ref{cor_main-spectral}(i) and proposition
\ref{prop_lim-qc-sober}(ii.b).
\end{proof}

\begin{proposition}\label{prop_T_0-criterion}
Let $(X,\cT)$ be a compact topological space, $\cU$ a family
of open and closed subsets of $X$, and endow $X$ with the
coarsest topology $\cT_\cU$ containing $\cU$. Then the
following conditions are equivalent :
\begin{enumerate}
\alphaenu
\item
$(X,\cT_\cU)$ is a $T_0$ topological space.
\item
$(X,\cT_\cU)$ is a spectral space.
\end{enumerate}
Moreover, if these conditions hold, then $(X,\cT_\cU)^c=(X,\cT)$.
\end{proposition}
\begin{proof} From remark \ref{rem_irred-comps}(vi) we
get (b)$\Rightarrow$(a).

(a)$\Rightarrow$(b): Since $\cT$ is compact, every $U\in\cU$
is compact as well in $\cT$, and therefore it is
quasi-compact in $\cT_\cU$. Thus, $(X,\cT_\cU)$ is
coherent, and it remains only to check that $X$ is
sober. To this aim, let $Z$ be any irreducible closed
subset of $(X,\cT_\cU)$, and set
$$
\cU_Z:=\{U\in\cU~|~U\cap Z\neq\emptyset\}
\qquad\text{and}\qquad
Z':=Z\cap\bigcap_{U\in\cU_Z}U.
$$
Since $\cT_\cU$ is $T_0$, the subset $Z'$ contains at most
one point of $X$, and to conclude, it suffices to check
that $Z'\neq\emptyset$. However, notice that $Z$ is closed
also in $\cT$, hence it is compact in this latter topology,
so the same holds for every finite intersection
$U_1\cap\cdots\cap U_k\cap Z$ with $U_1,\dots,U_k\in\cU_Z$.
We are then reduced to showing that every such finite
intersection is non-empty. Suppose that the latter fails,
and set $V_i:=Z\setminus U_i$ for every $i=1,\dots,k$;
then $Z=V_1\cup\cdots\cup V_k$, and since $Z$ is irreducible
in $\cT_\cU$, it follows that $Z=V_i$ for some $i\leq k$,
which is absurd.

Lastly, suppose that $(X,\cT_\cU)$ is spectral; we claim
that in this case the identity map of $X$ yields a
continuous map $i:(X,\cT)\to(X,\cT_\cU)^c$. Indeed, let
$T\subset X$ be a subset that is globally constructible
for the topology $\cT_\cU$; in view of lemma
\ref{lem_sorite-construct}(x.a,x.b), it suffices to show
that $T$ is open in the topology $\cT$.
Then lemma \ref{lem_sorite-construct}(ix.b) further reduces
to the case where $T=V\setminus V'$ for two quasi-compact
open subsets $V,V'$ of $(X,\cT_\cU)$. However, $V$ and $V'$
are finite unions of subsets of the form $U_1\cap\cdots\cap U_k$,
where $k\in\N$ is arbitrary, and $U_1,\dots,U_k\in\cU$;
it follows that $V$ and $V'$ are open and closed in $(X,\cT)$,
and therefore the same holds for their difference. Lastly,
since $(X,\cT)$ is separated and $(X,\cT_\cU)$ is quasi-compact
(theorem \ref{th_main-spectral}(i)), it follows that $i$ is a
homeomorphism, whence the proposition.
\end{proof}

\begin{definition}\label{def_special}
Let $X$ be a topological space, $x,x'\in X$ any two points.

\begin{enumerate}
\item
We say that {\em $x$ is a specialization of $x'$ in $X$}, or
equivalently, that {\em $x'$ is a generization of $x$ in $X$},
if $x$ lies in the topological closure of $\{x'\}$ in $X$.
\item
We say that
{\em $x$ is an immediate specialization of $x'$ in $X$},
or equivalently, that
{\em $x'$ is an immediate generization of $x$ in $X$}, if
$x\neq x'$, $x$ is a specialization of $x'$ in $X$, and any
other point of $X$ that is both a specialization of
$x'$ and a generization of $x$ in $X$ is equal to either
$x$ or $x'$.
\item
We denote by $X(x)$ the set of all generizations of
$x$ in $X$, and we endow $X(x)$ with the topology
induced by $X$.
\item
Let $f:X\to Y$ be any continuous map of topological
spaces. We say that $f$ is {\em specializing} (resp.
{\em generizing}) if the following holds. For every
$x\in X$ and every specialization (resp. generization)
$y$ of $f(x)$ in $Y$, there exists a specialization
(resp. generization) $x'$ of $x$ in $X$, such that
$f(x')=y$.
\end{enumerate}
\end{definition}

\begin{remark}\label{rem_specialize}
Let $X$ be a topological space, and $x$ any point of $X$. 

(i)\ \
Clearly $X(x)$ is the intersection of all open
neighborhoods of $x$ in $X$. Especially, if $X$
is coherent, $X(x)$ is a pro-constructible subset
of $X$. If $X$ is  $T_0$ ({\em i.e.} for any two
distinct points of $X$, at least one of them admits
an open neighborhood not containing the other), $x$
is the unique closed point of $X(x)$. If $X$ is
spectral, the same holds for $X(x)$, and the
inclusion map $j_x:X(x)\to X$ is quasi-compact
(corollary \ref{cor_procon-is-spec}); since $X(x)$
is quasi-separated, $j_x$ is also quasi-separated,
and even spectral, by remark \ref{rem_sorite-qcoh-maps}(iii).

(ii)\ \
Clearly, the relation
$$
x\leq x'
\quad\Leftrightarrow\quad
\text{$x$ is a specialization of $x'$ in $X$}
$$
defines a preordering on $X$ (see example
\ref{ex_universe}(iii)), and every continuous map
$f:X\to Y$ is also a morphism of preordered sets
$(X,\leq)\to(Y,\leq)$. If $X$ is $T_0$, then $(X,\leq)$
is even a partially ordered set. Suppose that $X$ is
$T_0$ and non-empty; in general, $(X,\leq)$ admits neither
maximal nor minimal points, but if $X$ is sober, it
follows easily from remark \ref{rem_irred-comps}(ii)
that $(X,\leq)$ admits maximal points. If $X\neq\emptyset$
is quasi-compact and $T_0$, then $(X,\leq)$ also admits
minimal elements: indeed, if $(x_i~|~i\in I)$ is any
(non-empty) totally ordered subset of $X$, let $\bar{\{x_i\}}$
be the topological closure of $\{x_i\}$ in $X$, for every
$i\in I$; then
$Z:=\bigcap_{i\in I}\bar{\{x_i\}}\neq\emptyset$, so any
point of $Z$ is smaller than every $x_i$, and then
the assertion follows as usual, by Zorn's lemma.

(iii)\ \
Suppose that $X$ is sober. Then, $X$ is noetherian
if and only if the following holds :
\begin{enumerate}
\alphaenu
\item
Every closed subset of $X$ has finitely many maximal
points.
\item
Every descending sequence $(x_n~|~n\in\N)$
in $(X,\leq)$ is stationary, {\em i.e.} it admits
$N\in\N$ such that $x_n=x_N$ for every $n\geq N$.
\end{enumerate}
Indeed, for a sober space, (a) and (b) are rephrasing
of the corresponding conditions of proposition
\ref{prop_noetherian}. It is also clear that, in
this situation, the dimension of $X$ is the supremum
of the lengths of all finite strictly descending sequences
in $(X,\leq)$ (see remark \ref{rem_irred-comps}(iv)).
\end{remark}

\begin{proposition}\label{prop_closed-under-spec}
Let $X$ be any locally spectral topological space,
$T\subset X$ a pro-constructible subset, and
$U\subset X$ an ind-constructible subset.
The following holds :
\begin{enumerate}
\item
The topological closure of\/ $T$ in $X$ is the set
of all specializations of all points of\/ $T$.
\item
$U$ is open in $X$ if and only if it contains the
generizations of all its points.
\item
$T$ is dense in $X$ if and only if it contains all
the maximal points of the partially ordered set
$(X,\leq)$ (notation of remark {\em\ref{rem_specialize}(ii)}).
\end{enumerate}
\end{proposition}
\begin{proof}(i): Denote by $\bar T$ the topological
closure of $T$ in $X$, and by $T^s$ the set of
all specializations in $X$ of points of $T$.
Clearly $T^s\subset\bar T$, so we need only prove
the converse inclusion. However, let $V$ be any
open subset of $X$; on the one hand, the topological
closure of $T\cap V$ in $V$ equals $V\cap\bar T$.
On the other hand, say that $x\in V$ is a specialization
in $X$ of a point $y\in T$; then clearly $y\in T\cap V$,
{\em i.e.} $x$ is a specialization in $V$ of $y$,
so $V\cap T^s$ is the set of specializations in $V$
of the points of $V\cap T$. Since $X$ admits a covering
consisting of spectral open subsets, we may then replace
$X$ by any of these subsets, and assume from start
that $X$ is spectral.

Now, suppose that there exists
$x\in \bar T\!\setminus\!T^s$, denote by $\cF$ the
family of all quasi-compact open neighborhoods of
$x$ in $X$, and notice that $\cF$ is cofiltered,
since $X$ is quasi-separated, and not empty, since
$X$ is quasi-compact. Since $x\in\bar T$, we must have
$W\cap T\neq\emptyset$ for every $W\in\cF$; on the
other hand, we have
$$
\bigcap_{W\in\cF}(W\cap T)=\emptyset
$$
since $x\notin T^s$. But both $T$ and the elements
of $\cF$ are closed in $X^c$, so $(W\cap T~|~W\in\cF)$
is a closed filter of $X^c$, therefore its center
cannot be empty (theorem \ref{th_main-spectral}(i)
and \eqref{subsec_filter}), a contradiction.

(ii) follows easily from (i), by considering the
pro-constructible subset $X\!\setminus\!U$.

(iii) is an immediate consequence of (i).
\end{proof}

\begin{example}\label{ex_boolean}
Let $X$ be a spectral space of dimension $\leq 0$. It follows
easily from proposition \ref{prop_closed-under-spec}(i)
that every quasi-compact open subset of $X$ is closed.
Then $X$ is separated : indeed, let $x,x'\in X$
be any two distinct points; since $X$ is sober, we may
assume that $x$ does not lie in the topological closure
of $\{x'\}$, in which case there exists a quasi-compact
open neighborhood $U$ of $x$ that does not contain $x'$,
and the foregoing implies that $X\setminus U$ is an open
neighborhood of $x'$, whence the contention. Conversely,
every separated and compact topological space whose
quasi-compact open subsets form a basis of the topology
is a spectral space of dimension $\leq 0$; indeed, such a
space is coherent, and it is obviously sober, since all
its points are closed.
A topological space with this properties is sometimes
called a {\em compact boolean space} : see \cite[p.168,169]{Kel},
where it is also explained that, by virtue of the Stone
representation theorem, such spaces are characterized as
those topological spaces which are homeomorphic to the
maximal spectra of (unital) boolean algebras. More
precisely, to such any compact boolean space $X$, one
attaches its {\em characteristic ring} $R_X$, which is
the boolean algebra of all continuous functions
$X\to\Z/2\Z$, and then $X$ is naturally homeomorphic to
$\Max\,R_X$.
\end{example}

\begin{corollary}\label{cor-pro-constr}
Let $f:X\to Y$ be a continuous map of locally
spectral topological spaces, and $T\subset X$
(resp. $S\subset Y$) any pro-constructible subset.
Then :
\begin{enumerate}
\item
If $f$ is quasi-compact, $f(T)$ is a
pro-constructible subset of\/ $Y$.
\item
If $f$ is spectral and generizing, the topological
closure in $X$ of $f^{-1}S$ is the preimage of the
topological closure of $S$ in $Y$.
\end{enumerate}
\end{corollary}
\begin{proof}(i): The assertion is local on $Y$, hence
we may assume that $Y$ is spectral; then $X$ admits
a finite covering consisting of open spectral subsets,
and we are easily reduced to the case where $X$ is
spectral as well. Now, by theorem
\ref{th_main-spectral}(i,ii), the subset $T$ is
closed, and therefore quasi-compact, in $X^c$.
Then $f(T)$ is quasi-compact in $Y^c$; hence it
is closed in $Y^c$, and therefore pro-constructible
in $Y$.

(ii): Let $\bar S$ denote the topological closure
of $S$ in $Y$; clearly the topological closure
of $f^{-1}S$ lies in $f^{-1}\bar S$, so we have
only to check the converse inclusion. Thus, say
that $x\in f^{-1}\bar S$, so $s:=f(x)\in\bar S$,
and therefore $s$ is a specialization of
a point $s'\in S$ (proposition
\ref{prop_closed-under-spec}(i)); since $f$
is generizing, it follows that there exists
a generization $x'$ of $x$ in $X$, such that
$f(x')=s'$. This shows that $x$ lies in the
set of specializations of the points of $f^{-1}S$,
and the contention follows.
\end{proof}

\begin{proposition}\label{prop_old-Groth}
Let $f:Y\to X$ be a surjective, spectral map of locally spectral
topological spaces, and suppose that $f$ is either specializing
or generizing. Then $X$ is the (topological) quotient of\/ $Y$,
under the equivalence relation induced by $f$.
\end{proposition}
\begin{proof} The assertion means that a subset $Z\subset X$ is
open (resp. closed) in $X$ if and only if $Z':=f^{-1}Z$ is closed
in $Y$. For any subset $Z$ of $X$ (resp. of $Y$), we
denote by $\overline Z$ the topological closure of $Z$ in $X$
(resp. in $Y$). Of course, we may assume that $Z'$ is closed
in $Y$, and we check that $Z$ is closed in $X$.
Now, $Z'$ is pro-constructible in $Y$ (lemma
\ref{lem_main-spectral}(ii)), so $f(Z')=Z$ is
pro-constructible in $X$ (corollary \ref{cor-pro-constr}(i)).
If $f$ is generizing, it follows that
$Z'=\bar{f^{-1}Z}=f^{-1}\bar Z$ (corollary
\ref{cor-pro-constr}(ii)); since $f$ is surjective, we deduce
that $Z=\bar Z$, {\em i.e.} $Z$ is closed, as required.
Lastly, suppose that $f$ is specializing, and notice that
$\bar Z$ is the set of all specializations of all points
of $Z$ (proposition \ref{prop_closed-under-spec}(i)); thus,
let $x\in Z$, $x'\in X$ a specialization of $x$ and $y\in Z'$
such that $f(y)=x$; since $f$ is specializing, there exists
$y'\in Y$ with $f(y')=x'$. But since $Z'$ is closed, we have
$y'\in Z'$, so $x'\in Z$, as required.
\end{proof}

\begin{corollary}\label{cor_min-max}
Let $X$ be a spectral topological space, whose set
$\Min\,X$ (resp. $\Max\,X$) of minimal (resp. maximal)
points is finite. If $Z\subset X$ is a subset closed
under specializations and generizations, then $Z$ is
open and closed.
\end{corollary}
\begin{proof} When $\Max\,X$ is finite, any such $Z$ is a
finite union of irreducible components, hence it is closed.
But the same applies to the complement of $Z$, whence the contention.
When $\Min\,X$ is finite, consider the continuous map
$$
g:Y:=\!\!\!\coprod_{x\in\Min\,X}\!\!\!X(x)\to X
$$
whose restriction to each $X(x)$ is the inclusion map.
Now, $Y$ is spectral, and $g$ is surjective, continuous
and spectral (remark \ref{rem_specialize}(i)), and it
is obviously also generizing; on the other hand, clearly
the assertion holds for $Y$, hence for $X$, by proposition
\ref{prop_old-Groth}.
\end{proof}

\begin{remark} Hochster defines in \cite{Hoch2}
an involution
$$
X\mapsto X^*
$$
of the category of spectral topological spaces
and spectral maps, with the following properties.
For every spectral space $X$, the set underlying
$X^*$ is the same as that underlying $X$, and the
constructible closed subsets of $X$ form a basis
of the topology of $X^*$. He shows that the partially
ordered set $(X^*,\leq)$ is the opposite of $(X,\leq)$
(\cite[Prop.8]{Hoch2}). Using this construction, it
is possible to give an alternative proof of corollary
\ref{cor_min-max}.
\end{remark}

\sset\subsubsection{}\label{subsec_not-pro-constr}
In the study of the valuation spectrum to be carried
out in section \ref{sec_Spv-of-ring} we shall encounter
certain subsets of a spectral space that are not
pro-constructible, but nevertheless are spectral
spaces, when endowed with the subspace topology.
In the following proposition we describe the general
principle underlying these constructions.

Let $X$ be any topological space, $S\subset X\times X$
a subset of {\em specializations of $X$} : {\em i.e.}
for every $(x,y)\in S$ the point $y$ is a specialization
in $X$ of the point $x$, in which case we shall also say
that $y$ is an {\em $S$-admissible specialization} of $x$,
and $x$ is an {\em $S$-admissible generization} of $y$.
We say that a subset $T\subset X$ is {\em $S$-closed}
if every $S$-admissible specialization of every point
of $T$ lies also in $T$. Furthermore, we shall say that
$S$ is a {\em transitive set of specializations}, if
the following holds : if $x,y,z\in X$ are any three
points, and $(x,y),(y,z)\in S$, then $(x,z)\in S$ as
well. We may then state :

\begin{proposition}\label{prop_Huber-crit}
Let $X$ be any spectral space, $S$ a given set of
specializations of $X$, set
$$
Y:=\{x\in X~|~
\text{$x$ has no proper $S$-admissible specializations}\}
$$
and endow $Y$ with the topology induced by the inclusion
map $Y\to X$. Suppose that :
\begin{itemize}
\item[(S1)]
Every $y\in Y$ has a fundamental system of open neighborhoods
in $X$ consisting of constructible $S$-closed open subsets.
\item[(S2)]
For every $x\in X$ there exists an $S$-admissible specialization
$y$ of $x$ that lies in $Y$.
\end{itemize}
Then the following holds :
\begin{enumerate}
\item
For every $x\in X$ there exists a unique point $r(x)\in Y$
that is an $S$-admissible specialization of $x$ in $X$.
\item
The topological space $Y$ is spectral, and the resulting
retraction
$$
r:X\to Y
$$
is spectral.
\item
More precisely, a subset $T\subset Y$ is constructible
in $Y$ if and only if\/ $r^{-1}T$ is constructible in $X$.
\end{enumerate}
\end{proposition}
\begin{proof}(i): Suppose that a given $x\in X$ has two
$S$-admissible specialization $y,y\in Y$; by remark
\ref{rem_irred-comps}(vi) we know that $X$ is $T_0$, so
we may assume that there exists an open neighborhood
$U$ of $y$ that does not contain $y'$. If $V\subset U$
is any open neighborhood of $y$, then $x\in V$, but
$V$ does not contain the $S$-admissible specialization
$y'$ of $x$, contradicting (S1).

(ii): Let $\cP$ (resp. $\cQ$) be the set of all open
subsets $T\subset Y$ such that $r^{-1}T$ is constructible
(resp. open and constructible) in $X$.

\begin{claim}\label{cl_Huber-crit}
$\cP$ is a basis of the topology of $Y$, and $Y$ is a
$T_0$ topological space.
\end{claim}
\begin{pfclaim} Since $X$ is $T_0$, condition (S1) implies
easily that the same holds for $Y$. Next, notice that we
have $U=r^{-1}r(U)$ for every open and $S$-closed subset
$U\subset X$, and therefore $U=r^{-1}(U\cap Y)$ for every
such $U$. Thus, $U\cap Y\in\cP$ for every open, constructible
and $S$-closed subset $U\subset X$, and the claim follows
from (S1).
\end{pfclaim}

Let now $\cT$ be the topology whose basis is $\cQ$;
clearly every element of $\cQ$ is open and closed in
the topology $\cT$. Moreover, since $X$ is quasi-compact,
it is easily seen that the same holds for $(Y,\cT)$.
Then proposition \ref{prop_T_0-criterion} and claim
\ref{cl_Huber-crit} show that $Y$ is spectral. It is
also clear that $r$ is continuous and quasi-compact, so
it is spectral, by remark \ref{rem_sorite-qcoh-maps}(iii).

(iii): Proposition \ref{prop_T_0-criterion} also shows
that $\cQ$ is the set of constructible subsets of $Y$,
whence the assertion.
\end{proof}

\begin{corollary}\label{cor_Huber-crit}
In the situation of proposition {\em\ref{prop_Huber-crit}},
suppose moreover that the following condition holds :
\begin{itemize}
\item[(S3)]
For every $x\in X$, any two $S$-admissible specializations
are {\em comparable}, {\em i.e.} if $y,z\in X$ are any
two $S$-admissible specializations of $x$, then either
$y$ is an $S$-admissible specialization of $z$ or else
$z$ is an $S$-admissible specialization of $y$.
\end{itemize}
Then, for every subset\/ $T\subset X$ the following
conditions are equivalent :
\begin{enumerate}
\alphaenu
\item
$T=r^{-1}r(T)$.
\item
$T$ is $S$-closed, and every $S$-admissible generization
of every point of\/ $T$ lies in $T$.
\end{enumerate}
\end{corollary}
\begin{proof} It is easily seen that (b)$\Rightarrow$(a).
Hence, suppose that (a) holds, let $x\in T$ be any point,
and $y\in X$ any $S$-admissible generization of $x$. By
(S3), either $x$ is an $S$-admissible specialization of
$r(y)$, or $r(y)$ is an $S$-admissible specialization of
$x$. If $r(y)$ is an $S$-admissible specialization of $x$,
then $r(x)=r(y)$, by proposition \ref{prop_Huber-crit}(i);
therefore $y\in T$, by (a). If $x$ is an $S$-admissible
specialization of $r(y)$, we must have $x=r(y)$, since
$r(y)\in Y$; then (a) yields again $y\in T$.

Lastly, let $y$ be an $S$-admissible specialization of
$x$. By (S3), either $r(x)$ is an $S$-admissible
specialization of $y$, or $y$ is an $S$-admissible
specialization of $r(x)$. If $r(x)$ is an $S$-admissible
specialization of $y$, then $r(x)=r(y)$, by proposition
\ref{prop_Huber-crit}(i); then (a) yields $y\in T$.
If $y$ is an $S$-admissible specialization of $r(x)$,
then $y=r(x)$, since $r(x)\in Y$; then again $y\in T$,
by (a).
\end{proof}

The foregoing notions shall be applied also to the study
of schemes and their underlying topological spaces,
which are always locally spectral (see remark
\ref{rem_Hochster}(ii)).
For future reference, we make the following :

\begin{definition}\label{def_coh-schemes}
Let $X$ be any scheme.
\begin{enumerate}
\item
We say that $X$ is {\em locally coherent}, if $\cO_X$
is a coherent sheaf of rings.
\item
We say that $X$ is {\em coherent}, if it is locally
coherent, quasi-compact and quasi-separated.
\end{enumerate}
\end{definition}

Notice that the topological space underlying a coherent
scheme is always spectral, by lemma
\ref{lem_charact-coh-sob-sp}(iv). We conclude this
section with a few miscellaneous results concerning
the topology of schemes.

\begin{proposition}\label{prop_max-is-qc}
Let $X$ be any reduced and coherent scheme, and endow the
set $\Max\,X$ of maximal points of $X$ with the topology
induced from $X$. Then $\Max\,X$ is quasi-compact.
\end{proposition}
\begin{proof} Clearly we have $\Max\,U=U\cap\Max\,X$ for
any open subset $U\subset X$; then, if $(U_i~|~i\in I)$
is any finite covering of $X$ consisting of affine open
subsets, it suffices to show that $\Max\,U_i$ is quasi-compact
for every $i\in I$. We may therefore assume from start that
$X$ is affine, say $X=\Spec\,A$ for some reduced coherent
ring $A$. We notice :

\begin{claim}\label{cl_youresovain}
Let $R$ be any reduced ring and $f\in R$ any element;
we have :
$$
\Max\,(\Spec\,R_f)=\Max\,(\Spec\,R)\cap\Spec\,R/\Ann_R(f).
$$
\end{claim}
\begin{pfclaim} Set $I:=\Ann_R(f)$ and let $\fp$ be any
minimal prime ideal of $R$; then $\fp\in\Max\,(\Spec\,R_f)$
if and only if $\fp\in\Max\,(\Spec\,R)$ and $f\notin\fp$.
Thus, it suffices to check that $f\notin\fp$ if and only
if $I\subset\fp$. However, if $f\notin\fp$, clearly
$I\subset\fp$. For the converse, notice that
$I_\fp=\Ann_{R_\fp}(f)$, and $R_\fp$ is a field, since $R$
is reduced; hence $I_\fp$ is either $R_\fp$ or $0$, depending
on whether the image of $f$ is invertible or not in $R_\fp$,
and the assertion follows easily.
\end{pfclaim}

Now, pick any covering $(U_i~|~i\in I)$ of $\Max\,X$
consisting of open subsets; we have to show that there
exists a finite subset $J\subset I$ such that $(U_i~|~i\in J)$
already covers $\Max\,X$. To this aim, we may assume
without loss of generality that every $U_i$ is of the
form $\Max\,(\Spec\,A_{f_i})$ for some $f_i\in A$.
Notice that $\Ann_A(f_i)$ is an ideal of finite type,
since $A$ is coherent; in view of claim \ref{cl_youresovain},
we deduce that for every $i\in I$ there exists a
constructible open subset $V_i$ of $X$ such that
$\Max\,X\!\setminus\!U_i=V_i\cap\Max\,X$. By assumption,
$\Max\,X\cap\bigcap_{i\in I}V_i=\emptyset$; since
every $V_i$ contains the generizations of all its
points, it follows that $\bigcap_{i\in I}V_i=\emptyset$.
By corollary \ref{cor_hair-split}(i) and theorem
\ref{th_main-spectral}, we deduce that there exists
a finite subset $J\subset I$ such that
$\bigcap_{i\in J}V_i=\emptyset$. Then clearly this
subset $J$ will do.
\end{proof}

\begin{lemma}\label{lem_flat-quot}
Let $A$ be a ring, $I\subset A$ an ideal. Then :
\begin{enumerate}
\item
The following are equivalent :
\begin{enumerate}
\item
The map $A\to A/I$ is flat.
\item
The map $A\to A/I$ is a localization.
\item
For every prime ideal $\fp\subset A$ containing $I$, we have $IA_\fp=0$,
especially, $V(I)$ is closed under generizations in $\Spec\,A$.
\end{enumerate}
\item
Suppose $I$ fulfills the conditions {\em (a)-(c)} of {\em(i)}.
Then the following are equivalent :
\begin{enumerate}
\item
$I$ is finitely generated.
\item
$I$ is generated by an idempotent.
\item
$V(I)\subset\Spec\,A$ is open.
\end{enumerate}
\end{enumerate}
\end{lemma}
\begin{proof}(i): Clearly (b)$\Rightarrow$(a). Also (a)$\Rightarrow$(c),
since every flat local homomorphism is faithfully flat.
Suppose that (c) holds; we show that the natural map $B:=(1+I)^{-1}A\to A/I$
is an isomorphism. Indeed, notice that $IB$ is contained in the
Jacobson radical of $B$, hence it vanishes, since it vanishes
locally at every maximal ideal of $B$.

(ii): Clearly (ii.c)$\Rightarrow$(ii.b)$\Rightarrow$(ii.a).
If (ii.a) holds, then $V(I)$ is constructible and closed
under generizations, by (ii.c), so it is open (proposition
\ref{prop_closed-under-spec}(ii)), {\em i.e.}
(ii.a)$\Rightarrow$(ii.c).
\end{proof}

\sset\subsubsection{}
For every ring $A$, we let $\cI(A)$ be the set of all ideals
$I\subset A$ fulfilling the equivalent conditions (i.a)-(i.c) of lemma
\ref{lem_flat-quot}, and $\cZ(A)$ the set of all closed subset
$Z\subset\Spec\,A$ that are closed under generizations in
$\Spec\,A$. In light of lemma \ref{lem_flat-quot}(i.c), we
have a natural mapping:
\set\begin{equation}\label{eq_from-I-to-Z}
\cI(A)\to\cZ(A)\quad:\quad I\mapsto V(I).
\end{equation}

\begin{lemma}\label{lem_from-I-to-Z}
The mapping \eqref{eq_from-I-to-Z} is a bijection,
whose inverse assigns to any $Z\in\cZ(A)$ the ideal $I[Z]$
consisting of all the elements $f\in A$ such that $fA_\fp=0$
for every $\fp\in Z$.
\end{lemma}
\begin{proof} Notice first that $I[V(I)]=I$ for every $I\in\cI(A)$;
indeed, clearly $I\subset I[V(I)]$, and if $f\in I[V(I)]$, then
$fA_\fp=0$ for every prime ideal $\fp$ containing $I$, hence the
image of $f$ in $A/I$ vanishes, so $f\in I$. To conclude the proof,
it remains to show that if $Z$ is closed and closed under generizations,
then $V(I[Z])=Z$. However, say that $Z=V(J)$, let $\fp\in Z$
be any prime ideal, and suppose that $f\in J$; then
$\Spec\,A_\fp\subset Z$, hence $f$ is nilpotent in $A_\fp$,
so there exist $n\in\N$ and an open neighborhood
$U\subset\Spec\,A$ of $\fp$ such that $f^n=0$ in $U$.
Since $Z$ is quasi-compact, finitely many such $U$
suffice to cover $Z$, hence we may find $n\in\N$ large
enough such that $f^n\in I[Z]$, whence the contention.
\end{proof}

The following result improves upon \cite[Cor.2.6]{Oh-Pe}.

\begin{proposition}\label{prop_ohm}
Let $R$ be a ring, $R'$ any finitely generated $R$-algebra,
$\phi:\Spec\,R'\to\Spec\,R$ the induced morphism of schemes,
$\Sigma\subset\Spec\,R$ any subset, and endow $\Sigma$ 
(resp. $\phi^{-1}\Sigma$) with the topology induced from
$\Spec\,R$ (resp. from $\Spec\,R'$). If\/ $\Sigma$ is a
noetherian topological space, then the same holds for
$\phi^{-1}\Sigma$.
\end{proposition}
\begin{proof} Set $X:=\Spec\,R$ and $X':=\Spec\,R'$;
we may find an integer $n\in\N$ and a closed immersion
$X'\to\Spec\,R[T_1,\dots,T_n]$, so it suffices to prove
the assertion for $R'=R[T_1,\dots,T_n]$, and by factoring
the map $R\to R'$ as a composition
$R\to R[T_1,\dots,T_{n-1}]\to R'$, an easy induction
reduces to the case where $R'=R[T]$.

Now suppose, by way of contradiction, that the
proposition fails, and let $\cF$ be the family
of all closed subsets $Z$ of $\Sigma$ such that
$\phi^{-1}Z$ is not noetherian. Since $\Sigma$ is
noetherian, and $\cF$ is not empty, $\cF$ admits
minimal elements. Let $Z$ be a minimal element of
$\cF$, and endow $Z$ with the topology induced from
$\Sigma$; then $Z$ is also noetherian (remark
\ref{rem_sorite-spectral}(iv)), so we may replace
$\Sigma$ by $Z$, and assume that no proper closed
subset of $\Sigma$ lies in $\cF$. We claim that in
this case, $\Sigma$ is irreducible. Indeed, let
$\Sigma_1,\dots,\Sigma_k$ be the finitely many
irreducible components of $\Sigma$ (proposition
\ref{prop_noetherian}); if $k>1$, we have $\Sigma_i\notin\cF$
for $i=1,\dots,k$, but then remark \ref{rem_sorite-spectral}(v)
easily implies that $\Sigma\notin\cF$ as well, a contradiction.
Hence, the topological closure $\bar\Sigma$ of $\Sigma$ in $X$
is irreducible in $X$ (lemma \ref{lem_irreducible}(i));
we endow $\bar\Sigma$ with the reduced closed subscheme
structure induced by the closed immersion $\bar\Sigma\to X$,
so $\bar\Sigma$ is isomorphic to the spectrum of some quotient
of $R$. We may then replace $X$ by $\bar\Sigma$, after which
we may assume that $R$ is a domain, and $\Sigma$ is dense in
$X$. In this case, denote by $K$ the field of fractions
of $R$. We notice :

\begin{claim}\label{cl_elementary-divs}
Let $A$ be a domain, $K$ the field of fractions of $A$,
$I\subset A[T]$ any ideal, and $p(T)\in I$ a monic
polynomial that generates $I\cdot K[T]$. Then $p(T)$
generates $I$.
\end{claim}
\begin{pfclaim} Consider any $g(T)\in I$, and let
$g(T)=q(T)\cdot p(T)+r(T)$ be the euclidean division
of $g$ by $p$ in $K[T]$. Since $p$ generates $I\cdot K[T]$,
we have $r=0$, and since $p$ is monic, it is easily seen
that $q\in A[T]$, whence the claim.
\end{pfclaim}

Now, suppose that $(Z_n~|~n\in\N)$ is a strictly descending
chain of closed subsets of $\phi^{-1}\Sigma$, and for each
$n\in\N$, let $\bar Z_n$ be the topological closure of
$Z_n$ in $\Spec\,R'$, and $I_n\subset R'$ the largest
ideal such $\Spec\,R'/I_n=\bar Z_n$. There follows
an increasing chain of ideals $(I_n~|~n\in\N)$.
The corresponding sequence $(I_n\cdot K[T]~|~n\in\N)$
is stationary and consists of principal ideals of
$K[T]$, so we may find $N\in\N$ and a monic polynomial
$p(T)$ such that $p(T)$ generates $I_n\cdot K[T]$
for every $n\geq N$. Moreover, we may find some
$f\in R\!\setminus\!\{0\}$ such that $p(T)\in R_f[T]$
(where $R_f:=R[f^{-1}]$), and we may even assume that
$p(T)\in I_n\cdot R_f[T]$ for every $n\geq N$, in
which case claim \ref{cl_elementary-divs} shows that
$p(T)$ generates $I_n\cdot R_f[T]$ for every $n\geq N$.
Set $U:=\Spec\,R_f$; we conclude that the chain of
subsets $(\bar Z_n\cap\phi^{-1}U~|~n\in\N)$ is stationary.
However,
$Z'_n:=Z_n\cap\phi^{-1}U=\bar Z_n\cap\phi^{-1}(U\cap\Sigma)$
for every $n\in\N$, so the chain $(Z'_n~|~n\in\N)$
is stationary as well. On the other hand, by construction
$U\neq\emptyset$, so $\Sigma':=\Sigma\!\setminus\!U$ is
a proper closed subset of $\Sigma$, therefore $\phi^{-1}\Sigma'$
is noetherian, and especially, the chain
$(Z_n\cap\phi^{-1}\Sigma'~|~n\in\N)$ is stationary.
Summing up, we deduce that the chain $(Z_n~|~n\in\N)$
is stationary, which is absurd.
\end{proof}

\subsection{Topological groups}
This section collects a few generalities on topological
groups that shall be used in later sections.

\begin{definition}
A {\em topological group} is a pair $(G,\cT_G)$ where
$(G,\cdot,e_G)$ is a group and $\cT_G$ a topology on
the set underyling $G$, such that :
\begin{itemize}
\item
the composition law of $G$ is a continuous map
$(G,\cT_G)\times(G,\cT_G)\to(G,\cT_G)$
\item
the rule $g\mapsto g^{-1}$ for every $g\in G$ yields
a continuous map $(G,\cT_G)\to(G,\cT_G)$.
\end{itemize}
A {\em morphism of topological groups}
$\phi:(G,\cT_G)\to(H,\cT_H)$ is a group homomorphism
$\phi:G\to H$ that is continuous for the topologies
$\cT_G$ and $\cT_H$.
\end{definition}

\begin{remark}\label{rem_top-groups-general}
Let $(G,\cT)$ be any topological abelian group.

(i)\ \
If $U\subset G$ is any open neighborhood of the neutral
element $0$, set $-U:=\{-g~|~g\in U\}$. Clearly $-U$ is
still an open neighborhood of $0$ in $G$, and thus the
same holds for $U\cap-U$. Especially, $G$ admits a
fundamental system of open neighborhoods of $0$ consisting
of subsets $U$ such that $U=-U$.

(ii)\ \
Let $\cU$ be any fundamental system of open neighborhoods
of $0$ in $G$. It is easily seen that a subset $V\subset G$
is open if and only if for every $g\in V$ there exists
$U\in\cU$ such that $g+U:=\{g+h~|~h\in U\}\subset V$.
Especially, $\cU$ determines completely the topology of $G$.
Notice also that, due to the continuity of the composition
law of $G$, for every $U\in\cU$ there exists $U'\in\cU$ such
that $U'+U':=\{g+g'~|~g,g'\in U'\}\subset U$.

(iii)\ \
Conversely, let $\cU$ be any family of subsets of $G$
such that :
\begin{enumerate}
\alphaenu
\item
$0\in U$ for every $U\in\cU$
\item
$U=-U$ for every $U\in\cU$
\item
for every $U\in\cU$ there exists $U'\in\cU$ such that
$U'+U'\subset U$
\item
for every $U,U'\in\cU$ there exists $U''\in\cU$ such
that $U''\subset U\cap U'$.
\item
for every $U\in\cU$ and every $g\in U$ there exists
$U'\in\cU$ such that $g+U'\subset U$.
\end{enumerate}
Then there exists a unique topology $\cT$ on $G$ such
that $\cU$ is a fundamental system of neighborhoods of
$0$ in $G$ relative to the topology $\cT$, and $(G,\cT)$
is a topological group. Indeed, the uniqueness is clear
from (ii), which also describes $\cT$ explicitly, and it
is easily seen that the composition law $G\times G\to G$
and the map $G\to G$ : $g\mapsto-g$ are both continuous
for this topology : details left to the reader.

(iv)\ \
Let $\phi:G'\to G$ be a group homomorphism with $G'$ also
abelian, and endow $G'$ with the topology $\cT'$ induced
by $\cT$. Then $(G',\cT')$ is a topological group. Indeed,
pick any fundamental system $\cU$ of open neighborhoods of
$0$ in $G$ fulfilling conditions (a)-(e) of (iii); it is
easily seen that the system
$\phi^{-1}\cU:=(\phi^{-1}U~|~U\in\cU)$ is a fundamental
system of open neighborhoods of the neutral element
$0'$ of $G'$ relative to the topology $\cT'$, and
clearly $\phi^{-1}\cU$ fulfills the same conditions
(a)-(e), whence the assertion.

(v)\ \
Likewise, let $\psi:G\to G''$ be a group homomorphism
with $G''$ also abelian; for any system $\cU$ of open
neighborhoods of $0$ as in (iv), set
$\psi(\cU):=(\psi(U)~|~U\in\cU)$; it is easily seen
that the system $\psi(\cU)$ fulfills again conditions
(a)-(e) of (iii), hence there exists a unique topology
$\cT''$ on $G''$ for which $\psi(\cU)$ is a fundamental
system of open neighborhoods of the neutral element
$0''\in G''$, and such that $(G'',\cT'')$ is a topological
abelian group. If $\psi$ is surjective, it is easily
seen that $\cT''$ agrees with the topology induced by
$\cT$ via $\psi$. More generally, $\psi(G)$ is an open
subgroup of $G''$ for the topology $\cT''$, and the
topology on $\psi(G'')$ induced by $\cT''$ agrees with
the quotient topology induced by $\cT$ via the surjection
$\psi:G\to\psi(G)$ : details left to the reader. Especially,
$\cT''$ is independent of the choice of $\cU$, and is the
finest topology on $G''$ such that
$\psi:(G,\cT)\to(G'',\cT'')$ is a morphism of topological
abelian groups.

(vi)\ \
Let $(G',\cT')$ be another topological abelian group,
and $\phi:G\to G'$ a group homomorphism. Then $\phi$
is continuous if and only if it is {\em continuous at
the point $0\in G$}, {\em i.e.} if and only if for every
open neighborhood $U'\subset G'$ of the neutral element
$0'\in G'$ there exists an open neighborhood $U\subset G$
of $0$ in $G$ with $\phi(U)\subset U'$. Indeed, the condition
is obviously necessary. Conversely, suppose this condition
holds, and let $V'\subset G'$ be any open subset; for
every $g\in\phi^{-1}V$ we may pick an open neighborhood
$U'$ of $0'$ in $G'$ such that $U'\subset-\phi(g)+V$,
and by assumption there exists an open neighborhood
$U\subset G$ of $0$ such that $U\subset\phi^{-1}(U')$,
whence $g+U\subset\phi^{-1}V$, so $\phi^{-1}V$ is open in
$G$, whence the contention.

(vii)\ \
Notice that $G$ is separated if and only if the subset
$\{0\}$ is closed in $G$. Indeed, the condition is
obviously necessary. For the converse, suppose that
$\{0\}$ is closed, and let $g\in G\setminus\{0\}$ be
any other point; it suffices to show that there exist
open neighborhoods $U$ and $U'$ respectively of $0$ and
of $g$, such that $U\cap U'=\emptyset$. However, by
assumption there exists an open neighborhood $V$ of
$0$ that does not contain $g$; then let $U$ be any
open neighborhood of $0$ such that $U=-U$ and
$U+U\subset V$, and set $U':=g+U$. Suppose $u\in U\cap U'$;
then there exists $u'\in U$ such that $u=g+u'$, whence
$g=u-u'\in U+U\subset V$, a contradiction.
\end{remark}

\begin{lemma}\label{lem_no-need-of-dense}
Let $G$ be a topological abelian group,
$I_0\subset G_0\subset G$ two subgroups, that we endow
with the topology induced from $G$. Also, let $\cT_0$
and $\cT$ be the quotient topologies of $\bar G_0:=G_0/I_0$
and respectively $\bar G:=G/I_0$. The following holds :
\begin{enumerate}
\item
$\cT_0$ agrees with the topology induced by $\cT$
via the inclusion map $j:\bar G_0\to\bar G$.
\item
Suppose that $I_0$ is closed in $G_0$, denote by
$I\subset G$ the topological closure of $I_0$ in $G$,
and endow $\bar G{}':=G/I$ with the quotient topology
$\cT'$. Then the induced map
$$
j':\bar G_0\to\bar G{}'
$$
is injective, and $\cT_0$ agrees with the topology induced
by $\cT'$ via $j'$.
\end{enumerate}
\end{lemma}
\begin{proof}(i): We remark, quite generally :

\begin{claim}\label{cl_trivial-open}
Let
$$
\xymatrix{ X' \ar[r]^-{g'} \ar[d]_{f'} & X \ar[d]^f \\
Y' \ar[r]^-g & Y
}$$
be a cartesian diagram of topological spaces. If $f$
is an open map, the same holds for $f'$.
\end{claim}
\begin{pfclaim} Every open subset of $X'$ is a union of
subsets of the form $U\times_YV=g'^{-1}U\cap f'^{-1}V$, for
arbitrary open subsets $U\subset X$ and $V\subset Y'$.
Then $f'(U\times_YV)=V\cap g^{-1}f(U)$ is open in $Y'$,
whence the claim.
\end{pfclaim}

We consider the diagram
$$
\xymatrix{ G_0 \ar[r] \ar[d]_{\pi_0} & G \ar[d]^\pi \\
\bar G_0 \ar[r]^-j & \bar G
}$$
whose vertical arrows are the projections, and whose top
horizontal arrow is the inclusion map. It is easily seen
that this diagram is cartesian, provided we endow $G_0$ and
$\bar G_0$ with the topologies induced by the inclusions
into $G$ and respectively $\bar G$. Then $\pi$ is an open
map, so the same holds for $\pi_0$, by claim
\ref{cl_trivial-open}, whence the assertion.

(ii): Since $I_0$ is closed in $G_0$, we have $I_0=G_0\cap I$,
so $j'$ is clearly injective. Let $\cT'_0$ be the topology on
$\bar G_O$ induced by $\cT'$ via $j'$; the identity map is a
continuous bijection
$$
(\bar G_0,\cT_0)\to(\bar G_0,\cT'_0)
$$
so it suffices to check that it is also a closed map.
Now, let $\pi':G\to\bar G{}'$ and $\pi'_0:G_0\to\bar G_0$
be the projections, and let $Z\subset\bar G_0$ be a subset
that is closed relative to the topology $\cT_0$; by definition,
this means that $\pi_0^{\prime -1}Z=Z+I_0$ is closed in $G_0$. Let
$(Z+I_0)^c$ be the topological closure of $Z+I_0$ in $G$;
it is easily seen that $(Z+I_0)^c=\pi^{\prime-1}\pi'((Z+I_0)^c)$,
so $\pi'((Z+I_0)^c)$ is closed in $\bar G{}'$, and hence
$\bar G_0\cap\pi'((Z+I_0)^c)$ is closed in $\bar G_0$,
relative to the topology $\cT'_0$; lastly
$$
\bar G_0\cap\pi'((Z+I_0)^c)=(G_0\cap(Z+I_0)^c)/I_0=(Z+I_0)/I_0=Z
$$
since $Z+I_0$ is closed in $G_0$. The assertion follows.
\end{proof}

\begin{lemma}\label{lem_descent-of-openness}
Let $f:G\to H$ and $g:H\to K$ be two continuous homomorphisms
of topological groups. If $g\circ f$ is an open map, the same
holds for $g$.
\end{lemma}
\begin{proof} Set $h:=g\circ f$, and let $U$ be any open
neighborhood of the neutral element $e_H$ in $H$; since
$h(f^{-1}U)\subset g(U)$, we see that $g(U)$ contains an
open neighborhood of the neutral element $e_K$ in $K$.
Now, let $x\in g(U)$ be any element, and pick $y\in U$
such that $h(y)=x$; let also $V$ be an open neighborhood
of $e_H$ in $H$ such that $y\cdot V\subset U$; it follows
that $g(y\cdot V)\subset g(U)$. On the other hand, $g(V)$
contains an open neighborhood $W$ of $e_K$, so $x\cdot W$
is an open neighborhood of $x$ contained in $g(U)$. Since
$x$ is arbitrary, the assertion follows.
\end{proof}

\begin{definition}\label{def_Cauchy-nets}
Let $G$ be any abelian topological group.

(i)\ \
A {\em Cauchy net} on $G$ is a family
$(g_\lambda~|~\lambda\in\Lambda)$ of elements of $G$ indexed
by a filtered partially ordered set $\Lambda$, such that
the following holds. For every open neighborhood $U$ of the
neutral element $0\in G$ there exists $\lambda\in\Lambda$
such that $g_\mu-g_\nu\in U$ for every $\mu,\nu\geq\lambda$.

(ii)\ \
We say that two Cauchy nets
$g_\bullet:=(g_\lambda~|~\lambda\in\Lambda)$
and $g'_\bullet:=(g'_{\lambda'}~|~\lambda'\in\Lambda')$ of $G$
are {\em equivalent}, and we write $g_\bullet\sim g'_\bullet$,
if the following holds. For every open neighborhood $U$ of
$0\in G$ there exist $\lambda\in\Lambda$ and
$\lambda'\in\Lambda'$ such that $g_\mu-g'_{\mu'}\in U$
for every $\mu\geq\lambda$ and every $\mu'\geq\lambda'$.

(iii)\ \
We say that an element $h\in G$ is a {\em limit point} for
the Cauchy net $g_\bullet:=(g_\lambda~|~\lambda\in\Lambda)$ if
$g_\bullet$ is equivalent to the trivial Cauchy net
$(h_{\lambda'}:=h~|~\lambda'\in\Lambda')$ indexed by a partially
ordered set $\Lambda'$ that has exactly one element. In this
case, we also say that $g_\bullet$ {\em converges} to $h$ in $G$.

(iv)\ \
We say that $G$ is {\em complete} if every Cauchy net
of $G$ admits a limit point. We denote by
$$
\mathsf{AbTopGrp}
\qquad\text{and}\qquad
\mathsf{AbTopGrp}_\mathrm{cp.sep.}
$$
respectively the category of abelian topological groups
(with morphisms given by the continuous group homomorphisms)
and its full subcategory whose objects are the complete
and separated topological groups.
\end{definition}

\begin{remark}\label{rem_general-complete}
(i)\ \
In view of remark \ref{rem_top-groups-general}(vii), it is
easily seen that an abelian topological group $G$ is separated
if and only if every Cauchy net of $G$ admits at most one
limit point. In this case, the unique limit point of a Cauchy
net $(g_\lambda~|~\lambda\in\Lambda)$ shall be often denoted
$$
\lim_{\lambda\in\Lambda}g_\lambda.
$$

(ii)\ \
Moreover, it is easily seen that the relation $\sim$
introduced by definition \ref{def_Cauchy-nets} is indeed
an equivalence relation on the set of all Cauchy nets on $G$.
Clearly, if $g_\bullet\sim g'_\bullet$, then an element $h\in G$
is a limit point for $g_\bullet$ if and only if it is a limit
point for $g'_\bullet$.

(iii)\ \
Furthermore, if $\phi:G\to G'$ is a continuous group
homomorphism of topological abelian groups, and
$(g_\lambda~|~\lambda\in\Lambda)$ a Cauchy net of $G$,
then $\phi(g_\bullet):=(\phi(g_\lambda)~|~\lambda\in\Lambda)$
is a Cauchy net of $G'$, and if $g'_\bullet$ is another
Cauchy net of $G$ with $g_\bullet\sim g'_\bullet$, then
$\phi(g_\bullet)\sim\phi(g'_\bullet)$.

(iv)\ \
Suppose that $G$ is complete and separated, and let $H\subset G$
be any subgroup, that we endow with the topology induced by $G$.
Then it is easily seen that $H$ is complete if and only
if it is topologically closed in $G$.
\end{remark}

\begin{theorem}\label{th_complete-top-grps}
{\em (i)}\ \ 
The inclusion functor
$\mathsf{AbTopGrp}_\mathrm{cp.sep}\to\mathsf{AbTopGrp}$
admits a left adjoint
$$
\mathsf{AbTopGrp}\to\mathsf{AbTopGrp}_\mathrm{cp.sep.}
\qquad
(G,\cT)\mapsto(G^\wedge,\cT^\wedge)
$$
called the {\em completion functor}.

{\em(ii)}\ \
For every topological abelian group $(G,\cT)$, the unit of
adjunction $i_G:G\to G^\wedge$ is called the {\em completion
map} of $G$. It has dense image, and $\cT$ agrees with the
topology induced by $\cT^\wedge$ via $i_G$.

{\em(iii)}\ \
Let $\phi:(G,\cT)\to(G',\cT')$ be a morphism of topological
abelian groups with dense image and with $G'$ complete and
separated, and such that $\cT$ agrees with the topology
induced by $\cT'$ via $\phi$. By the adjunction of {\em (i)},
the map $\phi$ factors through $i_G$ and a unique continuous
group homomorphism $\phi^\wedge:G^\wedge\to G'$, and $\phi^\wedge$
is an isomorphism of topological groups.
\end{theorem}
\begin{proof} Let $(G,+,0)$ be any topological abelian group.
We let $G^\wedge$ be the set of equivalence classes of
Cauchy nets of $G$. For any such Cauchy net $g_\bullet$,
we denote by $[g_\bullet]$ the class of $g_\bullet$ in
$G^\wedge$. We define a composition law on $G^\wedge$
as follows. Given two Cauchy nets
$g_\bullet:=(g_\lambda~|~\lambda\in\Lambda)$ and
$g'_\bullet:=(g'_{\lambda'}~|~\lambda'\in\Lambda')$, let
$\Lambda\times\Lambda'$ be the product of $\Lambda$
and $\Lambda'$ in the category of partially ordered
sets (whose underlying set is the cartesian product
of $\Lambda$ and $\Lambda'$); it is easily seen that
the family
$g''_{(\bullet,\bullet)}:=
(g_\lambda+g'_{\lambda'}~|~(\lambda,\lambda')\in\Lambda\times\Lambda')$
is also a Cauchy net of $G$, and we set
$$
[g_\bullet]+[g'_\bullet]:=[g''_{(\bullet,\bullet)}].
$$
One checks easily that $[g''_{(\bullet,\bullet)}]$ depends
only on the classes of $g_\bullet$ and $g'_\bullet$, so
we get by this rule a well defined mapping
$G^\wedge\times G^\wedge\to G^\wedge$. Clearly the neutral
element for this addition law is the class $[0]$ of any
constant Cauchy net $(0_\lambda~|~\lambda\in\Lambda)$
with $0_\lambda:=0$ for every $\lambda\in\Lambda$.
Also, set $-g_\bullet:=(-g_\lambda~|~\lambda\in\Lambda)$;
then it is easily seen that $[g_\bullet]+[-g_\bullet]=[0]$.
The associativity and commutativity of the addition law
are likewise immediate, so $G^\wedge$ is indeed an abelian
group.

\begin{claim}\label{cl_complete-ess-small}
The set $G^\wedge$ is essentially small.
\end{claim}
\begin{pfclaim} Let $g_\bullet:=(g_\lambda~|~\lambda\in\Lambda)$
be any Cauchy net of $G$. Fix any fundamental system $\cU$
of open neighborhoods of $0$ in $G$, and endow it with
the partial ordering such that $U\leq U'$ if and only
$U'\subset U$, for every $U,U'\in\cU$; for every $U\in\cU$
there exists $\lambda(U)\in\Lambda$ such that
$g_\lambda-g_\mu\in U$ for every $\lambda,\mu\geq\lambda(U)$.
It is easily seen that the family $(g_{\lambda(U)}~|~U\in\cU)$
is a Cauchy net of $G$ equivalent to $g_\bullet$. Thus, $G^\wedge$
has the cardinality of the set of equivalence classes
$G^\wedge_\cU$ of all Cauchy nets indexed by $\cU$; clearly
$\cU$ is small, so the same holds for $G^\wedge_\cU$, whence
the claim.
\end{pfclaim}

In light of claim \ref{cl_complete-ess-small}, we may replace
$G^\wedge$ by an isomorphic small group, that we shall still
denote by $G^\wedge$. We also have a natural group homomorphism
$$
i_G:G\to G^\wedge
\qquad
g\mapsto[g]
$$
that assigns to every $g\in G$ the class of any constant
Cauchy net with value $g$. Next, by remark
\ref{rem_top-groups-general}(i), we may find a fundamental
system $\cU$ of open neighborhoods of $0$ in $G$ such that
$-U=U$ for every $U\in\cU$; for every $U\in\cU$, we denote
by $U^\wedge\subset G^\wedge$ the subset consisting of the
classes of all Cauchy nets $(g_\lambda~|~\lambda\in\Lambda)$
fulfilling the following condition. There exists $U'\in\cU$
such that $g_\lambda+U'\subset U$ for every $\lambda\in\Lambda$.

\begin{claim}\label{cl_complete-U}
The system $\cU^\wedge:=(U^\wedge~|~U\in\cU)$ fulfills
conditions (a)-(e) of remark \ref{rem_top-groups-general}(iii).
\end{claim}
\begin{pfclaim} Indeed, (a) and (d) are obvious. To see
that (b) holds, let $(g_\lambda~|~\lambda\in\Lambda)$ be
any Cauchy net, and $U',U\in\cU$ such that
$g_\lambda+U'\subset U$ for every $\lambda\in\Lambda$;
then $-g_\lambda+(-U')\subset-U$, {\em i.e.}
$-g_\lambda+U'\subset U$ for every $\lambda\in\Lambda$,
whence the contention. Next, for any $U\in\cU$, pick
$U'\in\cU$ such that $U'+U'\subset U$; to check (c) it
suffices to show that $U'^\wedge+U'^\wedge \subset U^\wedge$.
Thus, say that $(g_\lambda~|~\lambda\in\Lambda)$ and
$(g'_{\lambda'}~|~\lambda'\in\Lambda')$ are two Cauchy
nets whose classes lie in $U'^\wedge$, and choose
$U''\in\cU$ such that $g_\lambda+U'',g'_{\lambda'}+U''\subset U'$
for every $\lambda\in\Lambda$ and every $\lambda'\in\Lambda'$;
it follows that $g_\lambda+g'_{\lambda'}+U''\subset U'+U'\subset U$
for every such $\lambda$ and $\lambda'$, as arrested.
To check (e), let $(g_\lambda~|~\lambda\in\Lambda)$ be
a Cauchy net, $U,U'\in\cU$ such that $g_\lambda+U'\subset U$
for every $\lambda\in\Lambda$, and
$(g'_{\lambda'}~|~\lambda'\in\Lambda)$ another Cauchy net
such that there exists $U''\in\cU$ with
$g'_{\lambda'}+U''\subset U'$ for every $\lambda'\in\Lambda'$;
it follows that $g_\lambda+g'_{\lambda'}+U''\subset U$ for
every such $\lambda$ and $\lambda'$, whence the contention.
\end{pfclaim}

From claim \ref{cl_complete-U} and remark
\ref{rem_top-groups-general}(iii) we deduce that there
exists a unique topology $\cT^\wedge$ on $G^\wedge$ for
which $\cU^\wedge$ is a fundamental system of open
neighborhoods of $0$, and $(G^\wedge,\cT^\wedge)$ is a
topological group. Moreover, we claim that
$$
i_G^{-1}(U^\wedge)=U
\qquad
\text{for every $U\in\cU$}.
$$
Indeed, the inclusion $U\subset i_G^{-1}(U^\wedge)$ is obvious.
For the converse, suppose that the class of the Cauchy net
$g_\bullet:=(g_\lambda~|~\lambda\in\Lambda)$ lies in $U^\wedge$,
and there exists $h\in G$ such that $i_G(h)=[g_\bullet]$.
This means that there exists $U'\in\cU$ such that
$g_\lambda+U'\subset U$ for every $\lambda\in\Lambda$ and
$\lambda_0\in\Lambda$ such that $h-g_\lambda\in U'$ for every
$\lambda\geq\lambda_0$; thus, $h\in g_\lambda+U'\subset U$,
whence the contention. Taking into account remark
\ref{rem_top-groups-general}(vi), it follows that $i_G$
is a continuous group homomorphism, and more precisely, the
topology of $G$ agrees with the topology induced by $\cT^\wedge$
via $i_G$. Furthermore, it is clear that the image of $i_G$
is dense in $G^\wedge$. Let us remark, more precisely :

\begin{claim}\label{cl_maximality-of-U-wedge}
(i)\ \
For every $U\in\cU$ and every open subset $V\subset G^\wedge$,
we have :
$$
V\subset U^\wedge
\qquad\Leftrightarrow\qquad
i^{-1}_GV\subset U.
$$

(ii)\ \
For every $U\in\cU$, the subset $U^\wedge$ lies in the topological
closure $i_G(U)^c$ of $i_G(U)$ in $G^\wedge$.
\end{claim}
\begin{pfclaim}(i): By the foregoing, we know already that
if $V\subset U^\wedge$, then $i^{-1}_GV\subset U$. Conversely,
suppose that $i^{-1}_GV\subset U$, and let
$h\in V\setminus U^\wedge$; then $h$ is the class of a Cauchy
net $h_\bullet:=(h_\lambda~|~\lambda\in\Lambda)$ with the
following property. For every $U'\in\cU$ there exist
$x_{U'}\in U'$ and $\lambda(U')\in\Lambda$ such that
$g_{U'}:=h_{\lambda(U')}+x_{U'}\notin U$. It is easily seen that
the family $g_\bullet:=(g_{U'}~|~U'\in\cU)$ is a Cauchy net of
$G$ equivalent to $h_\bullet$. We may thus replace $h_\bullet$
by $g_\bullet$, and assume from start that
$h_\lambda\in G\setminus U$ for every $\lambda\in\Lambda$.
Next, since $V$ is open in $G^\wedge$, and since we know
already that $i_G$ is continuous, there exists
$\lambda\in\Lambda$ such that $i_G(h_\lambda)\in V$; by
assumption, this implies that $h_\lambda\in U$, a contradiction.

(ii): Explicitly, $i_G(U)^c$ is the set of equivalence classes
of all Cauchy nets $g_\bullet:=(g_\lambda~|~\lambda\in\Lambda)$
such that $([g_\bullet]+V^\wedge)\cap i_G(U)\neq\emptyset$ for
every $V\in\cU$. Now, suppose that $[g_\bullet]\in U^\wedge$, so
we may assume that $g_\lambda\in U$ for every $\lambda\in U$,
and for any given $V\in\cU$ pick $V'\in\cU$ such that
$V'+V'\subset V$; we may then find $\lambda(V')\in\Lambda$
such that $g_\lambda-g_{\lambda'}\in V'$ for every
$\lambda,\lambda'\geq\lambda(V')$. Set
$\Lambda':=\{\lambda\in\Lambda~|~\lambda\geq\lambda(V')\}$,
$h:=g_{\lambda(V')}$ and $k_\lambda:=h-g_\lambda$ for every
$\lambda\in\Lambda'$. Clearly
$k_\bullet:=(k_\lambda~|~\lambda\in\Lambda')$ is a Cauchy net,
$h\in U$ and $[g_\bullet]+[k_\bullet]=i_G(h)$; lastly, we have
$k_\lambda+V'\subset V'+V'\subset V$ for every $\lambda\in\Lambda'$,
whence $[k_\bullet]\in V^\wedge$.
\end{pfclaim}

\begin{claim} $(G^\wedge,\cT^\wedge)$ is complete and
separated.
\end{claim}
\begin{pfclaim} Let $(g_{\bullet,i}~|~i\in I)$ a
family of Cauchy nets of $G$ indexed by a partially
ordered set $I$, such that $([g_{\bullet,i}]~|~i\in I)$
is a Cauchy net in $G^\wedge$. Hence, for every $i\in I$
we have a partially ordered set $\Lambda_i$ such that
$g_{\bullet,i}=(g_\lambda~|~\lambda\in\Lambda_i)$. For every
$U\in\cU$ we may find $i(U)\in I$ such that
$[g_{\bullet,i(U)}]-[g_{\bullet,j}]\in U^\wedge$ for every
$j\geq i(U)$. The latter means that there exist a Cauchy
net $(h_k~|~k\in K)$ of $G$ and $U'\in\cU$ such that
$h_k+U'\subset U$ for every $k\in K$, as well as
$\lambda(U)\in\Lambda_{i(U)}$, $\mu(U,j)\in\Lambda_j$ and
$k_0\in K$ such that $g_{\mu,j}-g_{\lambda,i(U)}-h_k\in U'$ for
every $k\geq k_0$, every $\mu\geq\mu(U,j)$ and every
$\lambda\geq\lambda(U)$. Therefore
$g_{\mu,j}\in g_{\lambda,i(U)}+h_k+U'\subset g_{\lambda,i(U)}+U$
for every such $j$, $\mu$ and $\lambda$. Lastly, we may also
assume that $g_{\lambda,i(U)}-g_{\lambda(U),i(U)}\in U$ for every
$\lambda\geq\lambda(U)$, whence
$g_{\mu,j}\in g_{\lambda(U),i(U)}+U+U$ for every $j\geq i(U)$
and every $\mu\geq\mu(U,j)$. Hence, the system
$l_\bullet:=(l_U:=g_{\lambda(U),i(U)}~|~U\in\cU)$ is a Cauchy net
in $G$ (for the partial ordering on $\cU$ such that $U\leq U'$
if and only if $U'\subset U$, for every $U,U'\in\cU$),
and it is easily seen that $[l_\bullet]$ is a limit point
of $([g_{\bullet,i}]~|~i\in I)$, which shows that $G^\wedge$
is complete.

Lastly, suppose that $g_\bullet:=(g_\lambda~|~\lambda\in\Lambda)$
is a Cauchy net such that $[g_\bullet]$ lies in $U^\wedge$ for
every $U\in\cU$. This means that for every such $U$ there
exists a Cauchy net
$h^U_\bullet:=(h^U_{\lambda'}~|~\lambda'\in\Lambda_U)$ such that
$h^U_{\lambda'}\in U$ for every $\lambda'\in\Lambda_U$ and such
that $g_\bullet\sim h^U_\bullet$. The latter condition implies
that there exist $\lambda(U)\in\Lambda$ and
$\lambda'(U)\in\Lambda_U$ such that $g_\lambda-h^U_{\lambda'}\in U$
for every $\lambda\geq\lambda(U)$ and every
$\lambda'\geq\lambda'(U)$. We conclude that
$g_\lambda\in U+U$ for every $\lambda\geq\lambda(U)$, so
$g_\lambda$ is equivalent to $[0]$, which shows that $G^\wedge$
is separated.
\end{pfclaim}

Next, consider a continuous group homomorphism $\phi:G\to H$
to a complete and separated topological abelian group $H$.
In view of remark \ref{rem_general-complete}(i,iii), it
follows easily that $\phi$ factors through $i_G$ and the group
homomorphism $\phi^\wedge:G^\wedge\to H$ that assigns to every
equivalence class $[g_\bullet]$ the unique limit point of the
Cauchy net $\phi(g_\bullet)$ in $H$. In view of remark
\ref{rem_top-groups-general}(vi), we need to check that
$\phi^\wedge$ is continuous at the point $0$. Thus, let
$V\subset H$ be any open neighborhood of the neutral point
$0_H$ of $H$, and let $g_\bullet:=(g_\lambda~|~\lambda\in\Lambda)$
be a Cauchy net in $G$ with $h:=\phi^\wedge([g_\bullet])\in V$;
we may find an open neighborhood $V'$ of $0_H$ in $H$ such
that $V'=-V'$ and $h+V'+V'+V'\subset V$, and set $U:=\phi^{-1}V'$,
so that $U=-U$. Now, let
$g'_\bullet:=(g'_{\lambda'}~|~\lambda'\in\Lambda')$ be any Cauchy
net in $G$ such that $[g'_\bullet]\in[g_\bullet]+U^\wedge$.
The condition means that there exist a Cauchy net
$(h_i~|~i\in I)$ in $G$ and an open neighborhood $U'$ of
$0$ in $G$ such that $h_i+U'\subset U$ for every $i\in I$,
and moreover there exist $\lambda_0\in\Lambda$,
$\lambda'_0\in\Lambda'$, and $i_0\in I$ such that
$g'_{\lambda'}-g_\lambda-h_i\in U'$ for every
$\lambda\geq\lambda_0$, every $\lambda'\geq\lambda'_0$
and every $i\geq i_0$. Thus, $g'_{\lambda'}\in g_\lambda+U$,
and consequently $\phi(g'_{\lambda'})\in\phi(g_\lambda)+V'$
for every $\lambda\geq\lambda_0$ and $\lambda'\geq\lambda'_0$.
Set $h':=\phi^\wedge([g'_\bullet])$; we may also assume that
$\phi(g_\lambda)-h,h'-\phi(g'_{\lambda'})\in V'$ for every
$\lambda\geq\lambda_0$ and every $\lambda'\geq\lambda'_0$.
We deduce that
$$
h'\in\phi(g'_{\lambda'})+V'\subset\phi(g_\lambda)+V'+V'
\subset h+V'+V'+V'\subset V
$$
{\em i.e.} $\phi^\wedge([g_\bullet]+U^\wedge)\subset V$, whence
the assertion. Lastly, since the image of $i_G$ is dense
in $G^\wedge$, the uniqueness of $\phi^\wedge$ is clear.

(iii): The existence of the continuous group homomorphism
$\phi^\wedge:G^\wedge\to G'$ has just been established. Next,
since $\phi$ has dense image, every $g'\in G'$ is the unique
limit point of a Cauchy net of $G'$ of the form
$(\phi(g_\lambda)~|~\lambda\in\Lambda)$ for a family
$g_\bullet:=(g_\lambda~|~\lambda\in\Lambda)$ of elements of
$G$, and since $\cT$ agrees with the topology induced by
$\cT'$, it is easily seen that $g_\bullet$ is a Cauchy net
in $G$. It then follows that $\phi^\wedge([g_\bullet])=g'$,
which shows that $\phi^\wedge$ is surjective. Let
$g_\bullet:=(g_\lambda~|~\lambda\in\Lambda)$ be a Cauchy
net in $G$ whose class $[g_\bullet]$ lies in the kernel
of $\phi^\wedge$; this means that
$(\phi(g_\lambda)~|~\lambda\in\Lambda)$ converges to $0$
in $G'$, and since $\cT$ is induced by $\cT'$, it follows
easily that $g_\bullet$ converges to $0$ in $G$, which
shows that $\phi^\wedge$ is also injective. It remains
therefore only to check that the topology on $G^\wedge$
induced by $\cT'$ via $\phi^\wedge$ is finer than $\cT^\wedge$.
To this aim, consider any $U\in\cU$; by assumption, there
exists an open neighborhood $W$ of $0$ in $G'$ such that
$\phi^{-1}W\subset U$. Set $V:=\phi^{\wedge\ -1}W$; then
$i_G^{-1}V\subset U$, and therefore $V\subset U^\wedge$,
by virtue of claim \ref{cl_maximality-of-U-wedge}(i). Since
$\cU^\wedge$ is a fundamental system of open neighborhoods
of $0$ in $G^\wedge$, the assertion follows.
\end{proof}

\begin{proposition}\label{prop_replaces-Mat-Th.8.1}
Let $(G,\cT_G)$ be a topological abelian group, $H\subset G$
a subgroup, and endow $H$ (resp. $G/H$) with the topology
$\cT_H$ (resp. $\cT_{G/H}$) induced by $\cT_G$ via the
inclusion map $j:H\to G$ (resp. via the projection
$\pi:G\to G/H$). Let also $i_G:G\to G^\wedge$ be the
completion map. Then we have :
\begin{enumerate}
\item
The resulting sequence of separated completions is exact :
$$
0\to H^\wedge\xrightarrow{\ j^\wedge\ }G^\wedge
\xrightarrow{\ \pi^\wedge\ }(G/H)^\wedge.
$$
\item
The topology $\cT^\wedge_H$ of $H^\wedge$ agrees with the one
induced by the topology $\cT^\wedge_G$ of $G^\wedge$ via $j^\wedge$.
\item
The image of $j^\wedge$ is the topological closure in $G^\wedge$
of $i_G(H)$.
\item
Endow $L:=G^\wedge/H^\wedge$ with the quotient topology $\cT_L$
induced by $G^\wedge$ via the projection $G^\wedge\to L$. Also,
let $k:L\to(G/H)^\wedge$ be the injective map deduced from
$\pi^\wedge$. Then $\cT_L$ agrees with the subspace topology
induced by $(G/H)^\wedge$ via $k$.
\item
Suppose that $G$ admits a countable fundamental system of
open neighborhoods of $0$. Then $\pi^\wedge$ is surjective.
\end{enumerate}
\end{proposition}
\begin{proof}(i): Let $g_\bullet:=(g_\lambda~|~\lambda\in\Lambda)$
be a Cauchy net of $H$; a direct inspection of the definitions
shows that $0$ is a limit point of $g_\bullet$ (relative to the
topology $\cT_H$) if and only if it is a limit point of
$i(g_\bullet)$, relative to the topology $\cT_G$, whence the
injectivity of $i^\wedge$. Next, by functoriality we have
$\pi^\wedge\circ j^\wedge=(\pi\circ j)^\wedge=0:H^\wedge\to(G/H)^\wedge$.
To conclude, suppose then that
$g_\bullet:=(g_\lambda~|~\lambda\in\Lambda)$ is a Cauchy net of
$G$ such that $\pi(g_\bullet):=(\pi(g_\lambda)~|~\lambda\in\Lambda)$
admits $0$ as limit point, and pick a fundamental system $\cU$
of open neighborhoods of $0$ in $G$; the assertion means that
for every $U\in\cU$ there exists $\lambda(U)\in\Lambda$ such
that $\pi(g_\lambda)\in\pi(U)$ for every $\lambda\geq\lambda(U)$,
{\em i.e.} $g_\lambda\in H+U:=\{h+u~|~h\in H,\ u\in U\}$ for
every such $\lambda$. Thus, for every $U\in\cU$ pick $h_U\in H$
such that $g_{\lambda(U)}-h_U\in U$; it follows that the families
$(g_{\lambda(U)}~|~U\in\cU)$ and $h_\bullet:=(h_U~|~U\in\cU)$ are
also Cauchy nets in $G$, and
$$
g_\bullet\sim(g_{\lambda(U)}~|~U\in\cU)\sim h_\bullet
$$
which shows that the class of $g_\bullet$ in $G^\wedge$ lies in
the image of $j^\wedge$, as stated.

(ii): Let $h_\bullet:=(h_\lambda~|~\lambda\in\Lambda)$ be a
Cauchy net of $H$ whose class $[h_\bullet]\in H^\wedge$ lies
in $j^{\wedge\ -1}(U^\wedge)$, for a given $U\in\cU$ (notation
of the proof of theorem \ref{th_complete-top-grps}); this
means that there exist a Cauchy net
$g_\bullet:=(g_{\lambda'}~|~\lambda'\in\Lambda')$ of $G$, and
$U'\in\cU$, such that $g_{\lambda'}+U'\subset U$ for every
$\lambda'\in\Lambda'$, and such that $h_\bullet\sim g_\bullet$.
Pick any $U''\in\cU$ such that $U''+U''\subset U'$; after
replacing $\Lambda$ and $\Lambda'$ by some cofinal subsets,
we may assume that $h_\lambda-g_{\lambda'}\in U''$ for every
$\lambda\in\Lambda$ and every $\lambda'\in\Lambda'$.
Whence, $h_\lambda+U''\in g_{\lambda'}+U''+U''\subset U$ for
every $\lambda\in\Lambda$, which shows that
$H^\wedge\cap U^\wedge\subset(H\cap U)^\wedge$. The converse
inclusion is clear, and the assertion follows.

(iii): On the one hand, the image of $j^\wedge$ is
a complete subgroup of $G^\wedge$, so it is topologically
closed (remark \ref{rem_general-complete}(iv)); on the
other hand, $i_G(H)$ is dense in in $H^\wedge$, whence
the contention.

(iv): Denote by $(\bar G,\bar\cT_G)$ the maximal separated
quotient of $G$ ({\em i.e.} the quotient $G/\{0\}^c$, where
$\{0\}^c$ denotes the topological closure of $\{0\}$ in
$G$), endowed with the quotient topology induced by $G$;
let also $\bar H$ be the topological closure of the
image of $H$ in $\bar G$, and endow $\bar H$ with the
topology $\bar\cT_{\!\!H}$ induced by $\bar\cT_{\!\!G}$ via
the inclusion map $\bar H\to\bar G$. Endow as well
$\bar G/\bar H$ with the quotient topology $\bar\cT_{\!\!G/H}$
induced by $\bar G$; we get a sequence of separated abelian
topological groups
$$
0\to\bar H\to\bar G\to\bar G/\bar H\to 0
$$
and after taking separated completions, a commutative
ladder of continuous group homomorphisms :
$$
\xymatrix{
0 \ar[r] & H^\wedge \ar[r] \ar[d] & G^\wedge \ar[r] \ar[d] &
(G/H)^\wedge \ar[d] \\
0 \ar[r] & \bar H{}^\wedge \ar[r] & \bar G{}^\wedge \ar[r] &
(\bar G/\bar H)^\wedge.
}$$ 
On the other hand, it is easily seen that $\cT_H$ (resp.
$\cT_{G/H}$) agrees with the topology induced by $\bar\cT_{\!\!H}$
(resp. by $\bar\cT_{\!\!G/H}$) via the natural map $H\to\bar H$
(resp. $G/H\to\bar G/\bar H$) : details left to the reader.
It follows immediately that the central and right-most
vertical arrows are isomorphisms of topological groups.
By the $5$-lemma, we deduce that the left-most vertical
arrow is a group isomorphism; but then it is also a
homeomorphism, by virtue of (ii). Summing up, we may
replace $G$ by $\bar G$ and $H$ by $\bar H$, and assume
from start that $G$ is separated and $H$ is topologically
closed in $G$; especially, the completion map $G\to G^\wedge$
is also injective. We now consider the commutative diagram
of continuous group homomorphisms
$$
\xymatrix{ G/H \ar[r]^-f \ar[rd]_{i_{G/H}} &
L \ar[r]^-{i_L} \ar[d]^k & L^\wedge \ar[ld]^g \\
& (G/H)^\wedge
}$$
where $i_{G/H}$ and $i_L$ are the completion maps.
By (iii) and lemma \ref{lem_no-need-of-dense}(ii), the
map $f$ is injective, and $\cT_{G/H}$ agrees with the
topology induced by $\cT_L$ via $f$. By theorem
\ref{th_complete-top-grps}(ii), the topology $\cT_L$
in turn is induced via $i_L$ by the topology $\cT^\wedge_L$
of its completion $L^\wedge$. Thus, $\cT_{G/H}$ is induced
by $\cT^\wedge_L$ via $i_L\circ f$; lastly, $i_L\circ f$
has dense image. Taking into account theorem
\ref{th_complete-top-grps}(iii), we conclude that $g$
is an isomorphism of topological groups, whence the
assertion.

(v): Let $\bar g_\bullet:=(\bar g_\lambda~|~\lambda\in\Lambda)$
be a Cauchy net in $G/H$, and pick a countable fundamental
system $\cU:=(U_n~|~n\in\N)$ of open neighborhoods of
$0$ in $G$. A simple induction shows that we may assume :
\set\begin{equation}\label{eq_good-fun-syst}
U_{n+1}+U_{n+1}\subset U_n
\qquad
\text{for every $n\in\N$}.
\end{equation}
Now, for every $n\in\N$ there exists $\lambda(n)\in\Lambda$
such that $\bar g_\lambda-\bar g_\mu\in\pi(U_n)$ for every
$\lambda,\mu\geq\lambda(n)$. It follows that the family
$(\bar g_{\lambda(n)}~|~n\in\N)$ is a Cauchy net in $G/H$
equivalent to $g_\bullet$ (for the standard ordering on $\N$);
thus, we may assume from start that the class $[\bar g_\bullet]$
in $(G/H)^\wedge$ is represented by a Cauchy {\em sequence}
$(\bar g_n~|~n\in\N)$ of $G/H$, such that
$\bar g_{n+1}-\bar g_n\in\pi(U_n)$ for every $n\in\N$.
Pick an arbitrary sequence $(g_n~|~n\in\N)$ in $G$ with
$\pi(g_n)=\bar g_n$ for every $n\in\N$; by construction,
for every $n\in\N$ there exists $h_n\in H$ such that
$g_{n+1}+h_n-g_n\in U_n$. Set $g'_n:=g_n+\sum_{i=0}^{n-1}h_i$
for every $n\in\N$; we then see that $g'_{n+1}-g'_n\in U_n$
for every $n\in\N$. Lastly, using \eqref{eq_good-fun-syst},
a simple induction shows more generally that
$g'_i-g'_j\in U_{n-1}$ whenever $i\geq j>0$. Therefore
$g'_\bullet:=(g'_n~|~n\in\N)$ is a Cauchy sequence in $G$,
and clearly $\pi^\wedge([g'_\bullet])=[g_\bullet]$, whence the
contention.
\end{proof}

\begin{corollary}\label{cor_product-and-completion}
Suppose that the topological group $G$ is the product of
a (small) family $(G_i~|~i\in I)$ of topological groups,
and for every $i\in I$ denote by $G^\wedge_i$ the separated
completion of $G_i$. Then the natural map
$j:G\to\prod_{i\in I}G^\wedge_i$ factors uniquely through the
completion map $G\to G^\wedge$ and an isomorphism of
topological groups :
$$
G^\wedge\isom\prod_{i\in I}G^\wedge_i.
$$
\end{corollary}
\begin{proof} It is easily seen that the product
$P:=\prod_{i\in I}G^\wedge_i$ is complete and separated;
also, $j$ has dense image, and the topology of $G$
is induced by that of $P$ via $j$. Then the assertion
follows from theorem \ref{th_complete-top-grps}(iii).
\end{proof}

\begin{corollary}\label{cor_limits-and-complete}
Let $G_\bullet:=(G_i~|~i\in\Ob(I))$ be any system of
topological abelian groups, indexed by any small
category $I$. We have :
\begin{enumerate}
\item
If $G_i$ is complete and separated for every $i\in\Ob(I)$,
then $\lim_IG_\bullet$ is a complete and separated topological
abelian group.
\item
Suppose that $I$ is cofiltered and countable (see example
{\em\ref{ex_count-filtered}}), and that the system $G_\bullet$
satisfies the following {\em Mittag-Leffler condition} :
\begin{itemize}
\item[(ML)]
For every $i\in\Ob(I)$ there exists a morphism $\phi:j\to i$ in
$I$ such that for every morphism $\psi:k\to j$ in $I$ we have
$$
\Img(G(\phi\circ\psi):G_k\to G_i)=\Img(G(\phi):G_j\to Gi).
$$
\end{itemize}
Then the natural map
$$
(\lim_IG_\bullet)^\wedge\to\lim_IG^\wedge_\bullet
$$
is an isomorphism of topological groups.
\end{enumerate}
\end{corollary}
\begin{proof}(i): Indeed, the limit $L$ of the system $G_\bullet$
is a closed subgroup of the product $P:=\prod_{i\in\Ob(I)}G_i$,
and the topology of $L$ is induced by that of $P$ via this
closed inclusion map; then the assertion follows from
corollary \ref{cor_product-and-completion} and remark
\ref{rem_general-complete}(iv).

(ii): In light of theorem \ref{th_complete-top-grps}(ii), it
is easily seen that the topology of $L$ is induced by that
of $L':=\lim_IG^\wedge_\bullet$, via the natural map
$\beta:L\to L'$. Taking into account (i) and theorem
\ref{th_complete-top-grps}(iii), we are then reduced to
checking that $\beta$ has dense image. Moreover, by virtue
of example \ref{ex_count-filtered}, we may assume that $I$
is the totally ordered set $\N^o$. Now, for every $i,j\in\N$
with $j\geq i$ let $\phi_{j,i}:G_j\to G_i$ be the transition
map; for every $i\in\N$ we set
$$
H_i:=\bigcap_{j\geq i}\Img\,\phi_{j,i}
$$
and we endow $H_i$ with the topology induced by $G_i$.
Condition (ML) easily implies that $\phi_{j,i}$ restricts to a
surjective continuous group homomorphism $H_j\to H_i$ for
every $i,j\in\N$ with $j\geq i$, and the natural morphisms
$$
\lim_{\N^o}H_\bullet\to L
\qquad
\lim_{\N^o}H^\wedge_\bullet\to L'
$$
are isomorphisms of topological groups. We may then replace
$G_\bullet$ by $H_\bullet$, and assume from start that the map
$\phi_{j,i}$ is a surjection, for every $j\geq i$.
In this case we come down to checking that the inclusion map
$H_i\to H^\wedge_i$ is dense for every $i\in\Ob(I)$, which we
know by theorem \ref{th_complete-top-grps}(ii).
\end{proof}

\begin{corollary}\label{cor_not-in-Bourbaki}
Let $G$ be any topological abelian group, $H\subset G$
an open subgroup, $H^\wedge$ and $G^\wedge$ the separated
completions of $H$ and $G$, and $i_G:G\to G^\wedge$ the
completion map. Then :
\begin{enumerate}
\item
$H^\wedge$ is an open subgroup of $G^\wedge$, and the natural
map $G/H\to G^\wedge/H^\wedge$ is a group isomorphism.
\item
$H=i_G^{-1}H^\wedge$.
\item
Suppose that $G$ admits a fundamental system
$(H_\lambda~|~\lambda\in\Lambda)$ of open neighborhoods of $0$,
consisting of open subgroups. Then
$(H^\wedge_\lambda~|~\lambda\in\Lambda)$ is a fundamental system
of open neighboroods of $0$ for $G^\wedge$.
\end{enumerate}
\end{corollary}
\begin{proof} Under these assumptions, $G$ induces the
discrete topology on the quotient $G/H$, hence the latter
is complete and separated; taking into account proposition
\ref{prop_replaces-Mat-Th.8.1}(i,ii), we get an exact
sequence
$$
0\to H^\wedge\to G^\wedge\to G/H
$$
from which both (i) and (ii) follow easily. Assertion (iii)
follows by inspecting the proof of theorem
\ref{th_complete-top-grps}(i).
\end{proof}

\begin{proposition}\label{prop_bilin-Bourbaki}
Let $G,H,K$ be three topological groups, and
$\phi:G\times H\to K$ a continuous and $\Z$-bilinear
map. Let also $i_G:G\to G^\wedge$, $i_H:H\to H^\wedge$ and
$i_K:K\to K^\wedge$ be the completion maps. Then there
exists a unique continuous $\Z$-bilinear map
$$
\phi^\wedge:G^\wedge\times H^\wedge\to K^\wedge
\qquad\text{such that}\qquad
\phi^\wedge\circ(i_G\times i_H)=i_K\circ\phi.
$$
\end{proposition}
\begin{proof} Let $x_\bullet:=(x_\lambda~|~\lambda\in\Lambda)$
and $y_\bullet:=(y_{\lambda'}~|~\lambda'\in\Lambda')$ be two
Cauchy nets in $G$ and respectively $H$, and let
$\Lambda\times\Lambda'$ be the product of $\Lambda$ and
$\Lambda'$ in the category of partially ordered sets;
{\em i.e.} the product of the underlying sets, endowed
with the partial ordering such that
$(\lambda,\lambda')\leq(\mu,\mu')$ if and only if
$\lambda\leq\mu$ and $\lambda'\leq\mu'$. Set
$x'_{(\lambda,\lambda')}:=x_\lambda$ and
$y'_{(\lambda,\lambda')}:=y_{\lambda'}$ for every
$(\lambda,\lambda')\in\Lambda\times\Lambda'$; it is easily
seen that the family $(x'_\mu~|~\mu\in\Lambda\times\Lambda')$
is a Cauchy net in $G$ equivalent to $x_\bullet$, and likewise
for the family $(y'_\mu~|~\mu\in\Lambda\times\Lambda')$.
Thus, we may replace $\Lambda$ and $\Lambda'$ by
$\Lambda\times\Lambda'$, after which we may assume that
$\Lambda=\Lambda'$.

\begin{claim}\label{cl_First-we-check}
The family $(\phi(x_\lambda,y_\lambda)~|~\lambda\in\Lambda)$
is a Cauchy net of $K$.
\end{claim}
\begin{pfclaim} Indeed, let $\cU_G$, $\cU_H$ and $\cU_K$ be
fundamental systems of open neighborhoods of the neutral
element for respectively $G$, $H$ and $K$, fulfilling
conditions (a)-(e) of remark \ref{rem_top-groups-general}(iii);
since $\phi$ is continuous for the product topology on
$G\times H$, for every $V\in\cU_K$ there exist $U_1\in\cU_G$
and $U_2\in\cU_H$ such that $\phi(U_1\times U_2)\subset V$.
On the other hand, there exists $\lambda(V)\in\Lambda$ such
that $x_\lambda-x_\mu\in U_1$ and $y_\lambda-y_\mu\in U_2$ for every
$\lambda,\mu\geq\lambda(V)$. It follows that
$$
v_{\lambda,\mu}^{(1)}:=\phi(x_\lambda-x_\mu,y_\lambda-y_{\lambda(V)}),
v_{\lambda,\mu}^{(2)}:=\phi(x_\mu-x_{\lambda(V)},y_\lambda-y_\mu)\in V
\qquad
\text{for every $\lambda,\mu\geq\lambda(V)$}.
$$
Moreover, the continuity of $\phi$ also yields
$\beta(V)\in\Lambda$ such that
$$
v^{(3)}_{\lambda,\mu}:=\phi(x_\lambda-x_\mu,y_{\lambda(V)}),
v^{(4)}_{\lambda,\mu}:=\phi(x_{\lambda(V)},y_\lambda-y_\mu)\in V
\qquad
\text{for every $\lambda,\mu\geq\beta(V)$}
$$
and clearly we may assume that $\beta(V)\geq\lambda(V)$
for every $V\in\cU_K$. We deduce 
$$
\phi(x_\lambda,y_\lambda)-\phi(x_\mu,y_\mu)=v^{(1)}_{\lambda,\mu}+
v^{(2)}_{\lambda,\mu}+v^{(3)}_{\lambda,\mu}+v^{(4)}_{\lambda,\mu}\in V+V+V+V
\qquad
\text{for every $\lambda,\mu\geq\beta(V)$}
$$
whence the contention.
\end{pfclaim}

In view of claim \ref{cl_First-we-check}, we denote by
$\phi(x_\bullet,y_\bullet):=
\lim_{\lambda\in\Lambda}\phi(x_\lambda,y_\lambda)\in K^\wedge$.
Next, we consider the group homomorphism
$$
\phi_{x_\bullet}:H\to K^\wedge
\qquad
y\mapsto\lim_{\lambda\in\Lambda}\phi(x_\lambda,y).
$$
It is easily seen that $\phi_{x_\bullet}$ is a continuous map,
hence it factors uniquely through $i_H$ and a continuous
group homomorphism $\phi^\wedge_{x_\bullet}:H^\wedge\to K^\wedge$.

\begin{claim}\label{cl_reduce-to-one-var}
$\phi(x_\bullet,y_\bullet)=
\phi^\wedge_{x_\bullet}(\lim_{\lambda\in\Lambda}y_\lambda)$.
\end{claim}
\begin{pfclaim} We resume the notation of the proof of
claim \ref{cl_First-we-check}, and notice that
$$
\phi(x_\lambda,y_\lambda)-\phi(x_\mu,y_\lambda)=
v^{(1)}_{\lambda,\mu}+v^{(3)}_{\lambda,\mu}\in V+V
\qquad
\text{for every $\lambda,\mu\geq\beta(V)$}.
$$
On the other hand, we may find $\gamma(V)\geq\beta(V)$
in $\Lambda$ such that
$$
\phi(x_\bullet,y_\bullet)-\phi(x_\lambda,y_\lambda)\in V
\qquad
\text{for every $\lambda\geq\gamma(V)$}
$$
as well as an element $\delta(V)\geq\gamma(V)$ in $\Lambda$
such that
$$
\phi(x_\mu,y_{\gamma(V)})-\phi^\wedge_{x_\bullet}(y_{\gamma(V)})\in V
\qquad
\text{for every $\mu\geq\delta(V)$}.
$$
Summing up, we conclude that
$$
\phi(x_\bullet,y_\bullet)-\phi^\wedge_{x_\bullet}(y_\mu)\in V+V+V+V
\qquad
\text{for every $\mu\geq\delta(V)$}
$$
and since $V$ is arbitrary, the assertion follows.
\end{pfclaim}

Claim \ref{cl_reduce-to-one-var} shows that
$\phi(x_\bullet,y_\bullet)$ depends only the class $y$ of
$y_\bullet$ in $H^\wedge$. Swapping the roles of $x_\bullet$
and $y_\bullet$, one may define a similar map
$\phi^\wedge_{y_\bullet}:G^\wedge\to K^\wedge$ such that
$\phi^\wedge_{y_\bullet}(x_\bullet)=\phi(x_\bullet,y_\bullet)$ as well,
which says likewise that $\phi(x_\bullet,y_\bullet)$ depends
only on the class $x$ of $x_\bullet$ in $G^\wedge$, so we
obtain a bilinear map
$\phi^\wedge:G^\wedge\times H^\wedge\to K^\wedge$ as sought, and
it remains only to check that $\phi^\wedge$ is continuous.
To this aim, define the fundamental systems $\cU_G^\wedge$,
$\cU_H^\wedge$, $\cU^\wedge_K$ of open neighborhoods of the
neutral element respectively in $G^\wedge$, $H^\wedge$ and
$K^\wedge$, as in the proof of theorem \ref{th_complete-top-grps};
we remark :

\begin{claim}\label{cl_many-wedges}
For every $W^\wedge\in\cU^\wedge_K$ there exist
$U^\wedge\in\cU^\wedge_G$, $V^\wedge\in\cU^\wedge_H$ with
$\phi^\wedge(U^\wedge\times V^\wedge)\subset W^\wedge$.
\end{claim}
\begin{pfclaim} Choose $W'\in\cU_H$ such that $W'+W'\subset W$;
then we may find $U\in\cU_G$, $V\in\cU_H$ such that
$\phi(U\times V)\subset W'$, whence $\phi(U\times V)+W'\subset W'$.
It follows easily that the corresponding open neighborhoods
$U^\wedge$, $V^\wedge$ will do.
\end{pfclaim}

Now, for a given $W^\wedge\in\cU^\wedge_K$, pick
$U^\wedge\in\cU^\wedge_G$, $V^\wedge\in\cU^\wedge_H$ as in
claim \ref{cl_many-wedges}; we may also assume that
$\phi^\wedge_{x_\bullet}(V^\wedge)\subset W^\wedge$ and
$\phi^\wedge_{y_\bullet}(U^\wedge)\subset W^\wedge$. It follows that
$$
\phi^\wedge((x+U^\wedge)\times(y+V^\wedge))-\phi^\wedge(x,y)\subset
\phi^\wedge(U^\wedge\times V^\wedge)+\phi^\wedge_{x_\bullet}(V^\wedge)+
\phi^\wedge_{y_\bullet}(U^\wedge)\subset W^\wedge+W^\wedge+W^\wedge
$$
whence the contention. Lastly, the uniqueness of $\phi^\wedge$
is clear, since the image of $G\times H$ is dense in
$G^\wedge\times H^\wedge$.
\end{proof}

\subsection{Topological rings and topological modules}
If $T$ is any topological space and $S\subset T$ any
subset, we shall denote by $S^c$ the topological closure
of $S$ in $T$.

\begin{definition}\label{def_top-ring}
(i)\ \
A {\em topological ring} is a pair $(A,\cT_A)$ where $A$
is a (commutative, unital) ring and $\cT_A$ a topology
on the set underlying $A$, such that :
\begin{itemize}
\item
the multiplication law of $A$ is a continuous map
$(A,\cT_A)\times(A,\cT_A)\to(A,\cT_A)$
\item
the pair $((A,+,0),\cT_A)$ is a topological group.
\end{itemize}

(ii)\ \
Let $(A,\cT_A)$ be a topological ring, $M$ an $A$-module,
$\cT_M$ a topology on the set underlying $M$. We say that
$(M,\cT_M)$ is a {\em topological $(A,\cT_A)$-module}
(briefly : a {\em topological $A$-module}) if
\begin{itemize}
\item
the scalar multiplication is a continuous map
$(A,\cT_A)\times(M,\cT_M)\to(M,\cT_M)$
\item
the pair $((M,+,0),\cT_M)$ is a topological abelian group.
\end{itemize}

(iii)\ \
Let $(A,\cT_A),(B,\cT_B)$ be two topological rings.
A {\em morphism of topological rings} is a ring homomorphism
$A\to B$ that is continuous for the topologies $\cT_A$ and
$\cT_B$. Such a morphism $(A,\cT_A)\to(B,\cT_B)$ is also
called a {\em topological $(A,\cT_A)$-algebra}. The topological
$A$-algebras form a category
$$
A\tdu\TopAlg
$$
whose morphisms are the continuous $A$-algebra homomorphisms.
Especially, the ring $\Z$ with its discrete topology is a
topological ring, and $\Z\tdu\TopAlg$ is the category of
topological rings and continuous ring homomorphisms.

(iv)\ \
Let $(A,\cT_A)$ be a topological ring, $B\to A$ a ring
homomorphism, $(M,\cT_M)$ a topological $A$-module.
We say that the topology $\cT_M$ is {\em $B$-linear}
if there exists a fundamental system of open neighborhoods of
$0\in M$ for the topology $\cT_M$, consisting of $B$-submodules
of $M$.

(v)\ \
Especially, in the situation of (iv) we say that the topology
$\cT_A$ is {\em $B$-linear}, if the same holds for the
topological $A$-module $(A,\cT_A)$. We say also simply that
$\cT_A$ is {\em linear} if it is $A$-linear. We denote by
$$
A\tdu\TopAlg_{B\tdu\mathrm{lin}}
$$
the full subcategory of $A\tdu\TopAlg$ whose objects are
the $B$-linear topological $A$-algebras.
\end{definition}

\begin{remark}\label{rem_I_bullet-adic}
(i)\ \ 
Let $A$ be a ring, and $\cT_A$ a topology on $A$ such that
$(A,\cT_A)$ is a topological group, for the addition law of
$A$. Let $\cU$ be a fundamental system of open neighborhoods
of $0$ in $A$; it is easily seen that $(A,\cT_A)$ is a
topological ring if and only if the system $\cU$ fulfills
the following conditions :
\begin{itemize}
\item[(a)]
For every $U\in\cU$ there exists $U'\in\cU$ such that
$U'\cdot U':=\{a\cdot b~|~a,b\in U'\}\subset U$.
\item[(b)]
For every $a\in A$ and every $U\in\cU$ there exists $U'\in\cU$
such that $a\cdot U'\subset U$.
\end{itemize}

(ii)\ \
Especially, if $I_\bullet:=(I_\lambda~|~\lambda\in\Lambda)$ is
any system of ideals of $A$ that is cofiltered by inclusion,
then there is a unique group topology $\cT_{I_\bullet}$ on $A$
for which $I_\bullet$ is a fundamental system of open neighborhoods
of $0$, and $(A,\cT_{I_\bullet})$ is a topological ring. We call
$\cT_{I_\bullet}$ the {\em $I_\bullet$-adic topology} on $A$.
Likewise, if $M$ is any $A$-module, the {\em $I_\bullet$-adic
topology} on $M$ is the unique group topology $\cT_M$ for which
$(I_\lambda M~|~\lambda\in\Lambda)$ is a fundamental system
of open neighborhoods of $0$ in $M$, and $(M,\cT_M)$ is
a topological $(A,\cT_{I_\bullet})$-module.
\end{remark}

\begin{remark}\label{rem_completion-of-topring}
Let $(A,\cT_A)$ be any topological ring and $(M,\cT_M)$
any topological $A$-module.

(i)\ \ 
Let $(A^\wedge,\cT^\wedge_A)$ be the separated completion
of the topological group underlying $A$; from proposition
\ref{prop_bilin-Bourbaki} we see that the multiplication
law of $A$ induces a continuous $\Z$-bilinear map
$A^\wedge\times A^\wedge\to A^\wedge$, which furnishes $A^\wedge$
with a natural structure of topological ring, such that the
completion map $A\to A^\wedge$ is a morphism of topological
rings. Likewise, the separated completion
$(M^\wedge,\cT_M^\wedge)$ of $(M,\cT_M)$ is naturally a
topological $A^\wedge$-module.

(ii)\ \
Suppose that $\cT_A$ is the linear topology defined by a
cofiltered system $I_\bullet:=(I_\lambda~|~\lambda\in\Lambda)$
of ideals of $A$, and let $I\subset A$ be any open ideal.
Also, suppose that $\cT_M$ is the $I_\bullet$-adic topology
on $M$. Then it follows easily from (i) and corollary
\ref{cor_not-in-Bourbaki}(iii), that $\cT^\wedge_M$ is the
$A^\wedge$-linear topology defined by the cofiltered system
$((I_\lambda M)^\wedge~|~\lambda\in\Lambda)$ of submodules of
$M^\wedge$. Moreover, $(IM)^\wedge$ is an open $A^\wedge$-submodule
in $M^\wedge$, the natural map $M/IM\to M^\wedge/(IM)^\wedge$
is an isomorphism, and the completion map $M\to M^\wedge$
identifies $(IM)^\wedge$ with the topological closure of
the image of $IM$ in $M^\wedge$ (corollary
\ref{cor_not-in-Bourbaki}(i)).

(iii)\ \
In the situation of (ii), suppose furthermore that for every
$\lambda\in\Lambda$ there exists $\mu\in\Lambda$ with
$I_\mu\subset I^2_\lambda$ (where $I^2_\lambda$ is the usual
square of the ideal $I_\lambda$). Then it is easily seen that
the $I_\bullet$-adic topology on $IM$ agrees with the topology
induced by $\cT_M$ on $IM$.

(iv)\ \
In the situation of (iii), suppose additionally that $\Lambda$
is countable and $I$ is finitely generated. Then the image of
$(IM)^\wedge$ in $M^\wedge$ equals $IM^\wedge$ : indeed, by assumption
we may find $k\in\N$ and a surjective $A$-linear map
$\phi:M^{\oplus k}\to IM$, and it is easily seen that $\phi$
is an open map; especially, the quotient topology induced
on $IM$ via the map $\phi$ agrees with $\cT_{IM}$. By
proposition \ref{prop_replaces-Mat-Th.8.1}(v), the
separated completion of $\phi$ is a surjective map 
$$
\phi^\wedge:(M^{\oplus k})^\wedge\to(IM)^\wedge
$$
and on the other hand we have a natural identification
$(M^{\oplus k})^\wedge\isom M^{\wedge\oplus k}$, whence the
contention (details left to the reader).
\end{remark}

\begin{example}\label{ex_clever-complete-mod}
(i)\ \
Let $A\neq 0$ be a ring, $\Sigma$ a set, and $\cF$ a {\em filter}
of $\Sigma$, {\em i.e.} a family of subsets of $\Sigma$ fulfilling
the following conditions : (a) for every $F,F'\in\cF$ we have
$F\cap F'\in\cF$; (b) for every $F\subset F'\subset\Sigma$ with
$F\in\cF$, we have $F'\in\cF$. Let us endow $A$ with the discrete
topology, and $M:=A^{(\Sigma)}$ with the linear topology $\cT_M$
defined by the system of submodules $(A^{(F)}~|~F\in\cF)$. Notice
that $\cT_M$ is separated if and only if
$\bigcap_{F\in\cF}F=\emptyset$. Notice as well the natural
isomorphism of $A$-modules
$$
M/A^{(F)}\isom A^{(\Sigma\setminus F)}
\qquad
\text{for every subset $F\subset\Sigma$}.
$$
For $F\subset F'\subset\Sigma$, these isomorphisms identify the
natural projection $M/A^{(F)}\to M/A^{(F')}$ with the $A$-linear
map $\pi_{F,F'}:A^{(\Sigma\setminus F)}\to A^{(\Sigma\setminus F')}$ :
$(a_\sigma~|~\sigma\in\Sigma\setminus F)\mapsto
(a_\sigma~|~\sigma\in\Sigma\setminus F')$. The completion of $M$
is then isomorphic to the inverse limit of the cofiltered
system $((A^{\Sigma\setminus F}~|~F\in\cF),\pi_{\bullet\bullet})$.
If $\cT_M$ is separated, a direct inspection shows that this
limit is represented by the submodule $L$ of $A^\Sigma$ consisting
of all sequences $(a_\sigma~|~\sigma\in\Sigma)$ with support
$\Lambda:=\{\sigma\in\Sigma~|~a_\sigma\neq 0\}$ satisfying the
following condition : $\Lambda\cap(\Sigma\setminus F)$ is a
finite set for every $F\in\cF$.

(ii)\ \
In the situation of (i), suppose moreover that for every
finite or countable subset $\Lambda\subset\Sigma$ there
exists $F\in\cF$ with $\Lambda\cap F=\emptyset$. Then $M$
is complete and separated. Indeed, the condition implies
that $\bigcap_{F\in\cF}F=\emptyset$, so $M$ is separated.
To check the completeness of $M$, say that
$a_\bullet:=(a_\sigma~|~\sigma\in\Sigma)\in L$; if the support
$\Lambda$ of $a_\bullet$ is not finite, we may find a countable
subset $\Lambda'\subset\Lambda$, and by assumption there
exists $F\in\cF$ such that $\Lambda'\cap F=\emptyset$.
But then $\Lambda'\subset\Lambda\cap(\Sigma\setminus\Lambda')$,
so $\Lambda\cap(\Sigma\setminus\Lambda')$ is not a finite
set, a contradiction. Thus $a_\bullet\in M$, whence the
assertion.
\end{example}

\begin{example}
Let $A$ be any topological ring whose topology is linear,
complete and separated, $I\subset A$ a closed ideal, and
endow $A/I$ with the quotient topology. In case $A$ admits
a fundamental system of open neighborhoods of $0$ consisting
of a {\em countable} family of ideals, then $A/I$ is still
complete and separated (proposition
\ref{prop_replaces-Mat-Th.8.1}(v)). But if this latter
condition is omitted, $A/I$ will still be separated, but
it will not necessarily be complete. For a counterexample,
take $A:=K[T_\lambda~|~\lambda\in\Lambda]$, where $K$ is a
non-zero ring, and $\Lambda$ is the (small) set of all
countable ordinal numbers, and let $I\subset A$ be the ideal
generated by $(T_\lambda-T_\mu~|~\lambda,\mu\in\Lambda)$. We
endow $A$ with the linear topology defined by the cofiltered
system of ideals $(J^n_\lambda~|~(\lambda,n)\in\Lambda\times\N)$,
where $J_\lambda\subset A$ is the ideal generated by
$(T_\mu~|~\mu\in\Lambda,\ \mu>\lambda)$, for every
$\lambda\in\Lambda$. Let $\Sigma\subset A$ be the subset of
monomials, which is in natural bijection with the set of all
maps $\phi:\Lambda\to\N$ such that $\phi^{-1}(\N\setminus\{0\})$
is a finite set. For every $(\lambda,n)\in\Lambda\times\N$,
set $\Sigma_{\lambda,n}:=J_\lambda^n\cap\Sigma$, and let $\cF$
be the set of all subsets $F\subset\Sigma$ such that
$\Sigma_{\lambda,n}\subset F$ for some
$(\lambda,n)\in\Lambda\times\N$. Then we have a natural
identification $A\isom K^{(\Sigma)}$, and the topology of
$A$ is defined by the cofiltered system of $K$-submodules
$(K^{(F)}~|~F\in\cF)$. It is easily seen that $\cF$ fulfills
the condition of example \ref{ex_clever-complete-mod}, hence
$A$ is complete and separated. Now, the quotient $A/I$ is
isomorphic to $K[T]$, endowed with its $T$-adic topology,
which is separated but not complete.
\end{example}

\begin{lemma}\label{lem_compl-and-int.clos}
Let $A$ be a topological ring, $B\subset A$ an open subring,
$C$ the integral closure of $B$ in $A$; endow $B$ and $C$
with the topologies induced by $A$, and denote by $A^\wedge$,
$B^\wedge$ and $C^\wedge$ the respective separated completions.
Then $C^\wedge$ is the integral closure of $B^\wedge$ in $A^\wedge$.
\end{lemma}
\begin{proof} Notice that $C^\wedge=B^\wedge\cdot C$, hence
$C^\wedge$ lies in the integral closure of $B^\wedge$ in $A^\wedge$,
and we are therefore reduced to checking that $B$ is integrally
closed in $A$ if and only if $B^\wedge$ is integrally closed in
$A^\wedge$. Let $j:A\to A^\wedge$ be the completion map. Suppose
first that $B$ is integrally closed in $A$, and let
$a\in A^\wedge$ be any element that is integral over $B^\wedge$;
pick any monic polynomial $P(T)\in B^\wedge[T]$ with $P(a)=0$.
Since $B$ is open in $A$, we may then find $a'\in A$ and a
monic polynomial $Q(T)\in B[T]$ such that
$$
j(a')-a,P(a'),j(Q(a'))-P(j(a'))\in B^\wedge.
$$
It follows that $Q(a')\in j^{-1}B^\wedge=B$ (corollary
\ref{cor_not-in-Bourbaki}(ii)). Set $R(T):=Q(T)-Q(a')$;
then $R(T)\in B[T]$ and $R(a')=0$, so $a'\in B$, by
assumption, and finally $a\in B^\wedge$, as required.
Conversely, if $a\in A$ is integral over $B$, it follows
easily that $j(a)$ is integral over $B^\wedge$, and
hence $j(a)\in B^\wedge$, if $B^\wedge$ is integrally
closed in $A^\wedge$; in this case, $a\in j^{-1}B^\wedge=B$,
which shows that $B$ is integrally closed in $A$.
\end{proof}

\sset\subsubsection{}\label{subsec_tensor-topol}
For any topological ring $A$, let $\cC_A$ denote the
full subcategory of $A\tdu\TopAlg$ whose objects are
the complete and separated topological $A$-algebras
whose topology is linear. Then the finite coproducts
of $\cC_A$ are representable. Indeed, consider more
generally any two topological $A$-modules $M$, $N$
whose topologies are $A$-linear; we denote by
$\cT^\otimes_{M,N}$ the $A$-linear topology on $M\otimes_AN$
defined by the system of submodules
$$
\Img(M'\otimes_AN+M\otimes_AN'\to M\otimes_AN)
$$
where $M'$ (resp. $N'$) ranges over the system of all
open $A$-submodules of $M$ (resp. of $N$), and following
\cite[Ch.0,\S7.7.2]{EGAI} we denote by
$$
M\,\hat\otimes_AN
$$
the separated completion of $(M\otimes_AN,\cT^\otimes_{M,N})$.
Notice that if $(M_\lambda~|~\lambda\in\Lambda)$ and
$(N_{\lambda'}~|~\lambda'\in\Lambda')$ are any two fundamental
systems of open submodules of $M$ and respectively $N$, we
have a natural isomorphism of topological $A$-modules
$$
M\,\hat\otimes_AN\isom\lim_{(\lambda,\lambda')\in\Lambda\times\Lambda'}
M/M_\lambda\otimes_AN/N_{\lambda'}
$$
where the tensor products $M/M_\lambda\otimes_AN/N_{\lambda'}$
are endowed with their discrete topologies.

By the same token, if $M'$ and $N'$ are two other topological
$A$-modules whose topologies are $A$-linear, and $f:M\to M'$
and $g:N\to N'$ are any two continuous $A$-linear maps, then
the map $f\otimes_Ag:M\otimes_AN\to M'\otimes_AN'$ is
continuous for the topologies $\cT^\otimes_{M,N}$ and
$\cT^\otimes_{M',N'}$, and therefore its separated completion
is a well defined continuous $A$-linear map
$$
f\,\hat\otimes_Ag:M\,\hat\otimes_AN\to M'\,\hat\otimes_AN'.
$$
Now, if $B$ and $C$ are any two topological $A$-algebras
whose topologies are linear, notice that the topology
$\cT^\otimes_{B,C}$ is defined by a cofiltered system of ideals
of $B\otimes_AC$, hence $(B\otimes_AC,\cT^\otimes_{B,C})$ is
a topological $A$-algebra.
We claim that if $B$ and $C$ are complete and separated,
then $B\,\hat\otimes_AC$ represents the coproduct of $B$
and $C$ in the category $\cC_A$. For the proof, consider
any two morphisms $f:B\to D$ and $g:C\to D$ in the category
$\cC_A$; there exists a unique morphism of $A$-algebras
$h:B\otimes_AC\to D$ such that $h(b\otimes c):=f(b)\cdot g(c)$
for every $b\in B$ and $c\in C$, and since $D$ is complete
and separated, it remains only to check that $h$ is
continuous for the topology $\cT^\otimes_{B,C}$. However,
let $I\subset D$ be any open ideal; then $J:=f^{-1}I$ (resp.
$K:=g^{-1}I$) is an open ideal of $B$ (resp. of $C$) and
clearly $h(J\otimes_AC+B\otimes_AK)\subset I$, whence
the contention.

\begin{definition}\label{def_bounded}
Let $(A,\cT)$ be any topological ring, and $S\subset A$
any subset.

(i)\ \
$S\subset A$ is {\em bounded} in $A$, if for every
open neighborhood $U$ of $0\in A$ there exists an
open neighborhood $V$ of $0\in A$ such that
$v\cdot s\in U$ for every $v\in V$ and $s\in S$.

(ii)\ \
For every $n\in\N$, let
$S(n):=\{a_1\cdots a_n~|~a_1,\dots,a_n\in S\}$.
We say that $S$ is {\em power bounded} in $A$ if
the subset $\bigcup_{n\in\N}S(n)$ is bounded in $A$.
An element $a\in A$ is {\em power bounded} (resp.
{\em topologically nilpotent}) in $A$, if the subset
$\{a\}$ is power bounded in $A$ (resp. if the
sequence $(a^n~|~n\in\N)$ converges to $0$ in the
topology $\cT$). We denote by
$$
A^\circ
\qquad\text{and}\qquad
A^{\circ\circ}
$$
respectively the subset of all power-bounded elements
and the subset of all topologically nilpotent elements of $A$.
An ideal $I\subset A$ is {\em topologically nilpotent}, if
for every open neighborhood $U$ of $0$ in $A$ there
exists $n\in\N$ such that $I^n\subset U$.

(iii)\ \
We say that $(A,\cT)$ is {\em adic} (resp. {\em c-adic})
if there exists an ideal $I$ of $A$ such that
$$
(I^n~|~n\in\N)
\qquad
\text{(resp.\ \ $((I^n)^c~|~n\in\N)$\ )}
$$
is a fundamental system of open neighborhood of zero in $A$.
Then, we call any such $I$ an {\em ideal of adic definition}
(resp. {\em an ideal of c-adic definition}) of $(A,\cT)$,
and we also say that $\cT$ is {\em $I$-adic} (resp.
{\em $I$-c-adic}). If $M$ is any $A$-module, the
{\em $I$-adic topology} on $M$ is the (unique) $A$-linear
topology such that $(I^nM~|~n\in\N)$ is a fundamental system
of open neighborhoods of $0$.

(iv)\ \
We say that $(A,\cT)$ is {\em f-adic} if there exists
an open subring $A_0\subset A$ that is adic for the
topology induced by $\cT$, and such that $A_0$ admits
a finitely generated ideal of adic definition.

(v)\ \
Let $(A,\cT)$ be any f-adic ring. An open subring $A_0$
of $A$ is called a {\em ring of definition of $A$}, if
$\cT$ induces a linear topology on $A_0$.

(vi)\ \
We say that $(A,\cT)$ is a {\em Tate ring} if it is
f-adic and $A^\times\cap A^{\circ\circ}\neq\emptyset$.
\end{definition}

\begin{remark}\label{rem_something-on-bdd}
Let $(A,\cT)$ be any topological ring.

(i)\ \ 
Let $S_1,S_2\subset A$ be two bounded subsets. Clearly
$S_1\cup S_2$ is bounded. Moreover, also
$S_1S_2:=\{s_1s_2~|~s_1\in S_1,\ s_2\in S_2\}$ is
bounded. Indeed, for any given open neighborhood $U$ of
$0$ in $A$, let us pick an open neighborhood $V$ of $0$
in $A$ such that $S_1V\subset U$, and an open neighborhood
$V'$ of $0$ in $A$ such that $S_2V'\subset V$; then
$S_1S_2V'\subset U$, whence the contention.

(ii)\ \
Let $T_1,\dots,T_k\subset A$ be any finite family of subsets.
Then $T:=T_1\cup\cdots\cup T_k$ is power bounded in $A$ if
and only if every $T_1,\dots,T_k$ is power bounded in $A$.
Indeed, set $S_i:=\bigcup_{n\in\N}T_i(n)$ for $i=1,\dots,k$
(notation of definition \ref{def_bounded}(ii)); by definition,
$T$ is power bounded if and only if the product $S_1\cdots S_k$
is bounded, so the assertion follows from (i).

(iii)\ \
Let $T\subset A$ be any subset. Then $T$ is bounded
(resp. power bounded) in $A$ if and only if the same
holds for the topological closure $T^c$ of $T$ in $A$.
Indeed, suppose that $T$ is bounded, and let $U\subset A$
be any open neighborhood of $0\in A$; pick an open
neighborhood $U'$ of $0$ in $A$ such that
$U'+U'\subset U$. By assumption, there exists an open
neighborhood $V$ of $0$ in $A$ such that
$T\cdot V\subset U'$, whence $T^c\cdot V\subset U'^c$;
but notice that $U'^c\subset U'+U'$ : indeed, if $x\in U'^c$,
we have $(x-U')\cap U'\neq\emptyset$, whence $x\in U'+U'$.
Summing up, we get $T^c\cdot V\subset U$, which shows that
$T^c$ is bounded. Next, suppose that $T$ is power bounded,
{\em i.e.} $T':=\bigcup_{n\in\N}T(n)$ is bounded, and notice
that $\bigcup_{n\in\N}T^c(n)\subset T'^c$; by the foregoing
$T'^c$ is bounded, so $T^c$ is power bounded.

(iv)\ \
Let $a\in A$ be any element, and $n\geq 1$ any integer.
Then $a$ is topologically nilpotent if and only if the
same holds for $a^n$. Indeed, obviously if $a$ is topologically
nilpotent, the same holds for $a^n$. Conversely, suppose
that $a^n$ is topologically nipotent; since the map
$A\to A$ : $b\mapsto a^rb$ is continuous for every $r\in\N$,
it follows that the sequences $(a^{kn+r}~|~k\in\N)$ converge
to $0$ in $A$ for every $r=0,\dots,n-1$; the claim is an
immediate consequence.

(v)\ \
Let $a,b\in A$ be any two elements. If $a$ is power bounded
and $b$ is topologically nilpotent, $ab$ is topologically
nilpotent. Indeed, for any open neighborhood $U$ of $0$ in
$A$ there exists an open neighborhood $V$ of $0$ in $A$ such
that $\{b^k~|~k\in\N\}\cdot V\subset U$, and on the
other hand, there exists $n\in\N$ such that $a^k\in V$ for
every integer $k\geq n$; thus, $(ab)^k\in U$ for every
$k\geq n$, whence the claim.
\end{remark}

\begin{remark}\label{rem_someth-on-bdd-in-Z-lin}
Let $(A,\cT)$ be a topological ring whose topology $\cT$
is $\Z$-linear.

(i)\ \
The subset $A^{\circ\circ}$ is an additive subgroup of $A$.
Indeed, let $a,b\in A^{\circ\circ}$ be any two elements, $U$
any open additive subgroup in $A$, and pick an open
neighborhood $V$ of $0$ in $A$ such that $V\cdot V\subset U$;
we may find $r\in\N$ such that $a^i,b^i\in V$ for every
$i\geq r$. Then, we may also find an integer $s\geq r$
such that $a^ib^j\in U$ for every $i,j\in\N$ such that
$i+j\geq s$ and $\min(i,j)<r$. It follows easily that
$(a+b)^n\in U$ for every $n\geq s$, and since $U$ is
arbitrary, we conclude that $a+b$ is topologically
nilpotent.

(ii)\ \
For any subset $T\subset A$, denote by $\La T\Ra\subset A$
the additive subgroup generated by $T$ in $A$. Then $T$
is bounded (resp. power bounded) if and only if the same
holds for $\La T\Ra$. Indeed, suppose that $T$ is bounded,
and let $U\subset A$ be any open additive subgroup; by
assumption there exists an open neighborhood $V$ of $0$
in $A$ such that $T\cdot V\subset U$. But then clearly
$\La T\Ra\cdot V\subset U$, so $\La T\Ra$ is bounded.
Likewise, suppose that $T$ is power bounded, and set
$S:=\bigcup_{n\in\N}T(n)$ (notation of definition
\ref{def_bounded}(ii)); since $S$ is bounded, we know
already that the same holds for $\La S\Ra$. However,
clearly $\La S\Ra=\bigcup_{n\in\N}\La T\Ra(n)$, so $\La T\Ra$
is power bounded.

(iii)\ \
Let $T\subset A$ be any power bounded subset; in light
of (ii), it is easily seen that the $\Z$-subalgebra
$\Z[T]\subset A$ generated by $T$ is bounded in $A$,
and especially, every element of $\Z[T]$ is power bounded.
Combining with remark \ref{rem_something-on-bdd}(ii), we
deduce that $A^\circ$ is the filtered union of all the bounded
subrings of $A$. Especially, $A^\circ$ is a subring of $A$.

(iv)\ \
Moreover, $A^\circ$ is integrally closed in $A$, and
$A^{\circ\circ}$ is a radical ideal of $A^\circ$. Indeed,
say that $x\in A$ is integral over $A^\circ$, so that
there exist $a_1,\dots,a_n\in A^\circ$ (for some $n\in\N$)
such that $x^n+a_1x^{n-1}+\cdots+a_n=0$, and set
$B:=\Z[a_1,\dots,a_n]\subset A$; by (iii), the subring
$B$ is bounded in $A$, and it is easily seen that
$B[x]=\sum^{n-1}_{i=0}Bx^i$, so $B[x]$ is bounded as well,
therefore $B[x]\subset A^\circ$, which shows the first
assertion. Lastly, (i) and remark
\ref{rem_something-on-bdd}(v) imply easily that
$A^{\circ\circ}$ is an ideal of $A^\circ$, and remark
\ref{rem_something-on-bdd}(iv) shows that $A^{\circ\circ}$
is a radical ideal.

(v)\ \
Lastly, if the topology $\cT$ is also complete and
separated, $A^{\circ\circ}$ lies in the Jacobson radical
of $A^\circ$. Indeed, in this case, for every
$a\in A^{\circ\circ}$ the series $1-a+a^2-a^3+\cdots$
converges in $A$ to a unique element $b$ such that
$b\cdot(1+a)=1$. Now, $b$ lies in the topological
closure $\Z[a]^c$ of $\Z[a]\subset A$, and notice
that $\Z[a]$ is the additive subgroup generated by
the subset $T:=\{a^n~|~n\in\N\}$; since $a$ is power
bounded, obviously $T$ is a power bounded subset,
and then the same holds for $\Z[a]^c$, by (ii) and
remark \ref{rem_something-on-bdd}(iii). Thus, $b\in A^\circ$,
whence the contention.
\end{remark}

\begin{lemma}\label{lem_adic-on-cartesian}
Let $A$ be a ring, and $I\subset A$ an ideal; consider
a cartesian diagram of topological $A$-modules
$$
\xymatrix{ M_0 \ar[r] \ar[d] & M_1 \ar[d]^{f_1} \\
M_2 \ar[r]^-{f_2} & M_3
}$$
and suppose that the topology of $M_3$ is discrete, and the
topologies of $M_1$ and $M_2$ are $I$-adic. Then the topology
of $M_0$ is $I$-adic as well.
\end{lemma}
\begin{proof} By assumption, we may find $c\in\N$ such
that $f_1(I^cM_1)=f_2(I^cM_2)=0$, and it follows easily that
$I^nM_1\times I^nM_2\subset M_0$ for every $n\geq c$.
Especially, $(I^nM_1\times I^nM_2~|~n\geq c)$ is a fundamental
system of open neighborhoods of $0$ in $M_0$. We also see that
$$
I^{n+c}M_1\times I^{n+c}M_2=I^n(I^cM_1\times I^cM_2)
\subset I^nM_0\subset I^nM_1\times I^nM_2
\qquad
\text{for every $n\in\N$}
$$
whence the claim.
\end{proof}

\begin{lemma}\label{lem_fontaine}
Let $A$ be a topological ring, $(M,\cT_M)$ a topological
$A$-module whose topology is $A$-linear, complete
and separated, and $I\subset A$ an ideal such that the
$I$-adic topology on $M$ is finer than $\cT_M$. Suppose
moreover that either one of the following conditions holds :
\begin{enumerate}
\alphaenu
\item
$I^nM$ is closed in the topology $\cT_M$, for every $n\in\N$.
\item
$I$ is finitely generated.
\end{enumerate}
Then $M$ is separated and complete for the $I$-adic topology.
\end{lemma}
\begin{proof} Let $(J_\lambda~|~\lambda\in\Lambda)$ be a
fundamental system of open submodules in $M$; by assumption,
for every $\lambda\in\Lambda$ there exists $n(\lambda)\in\N$
such that $I^{n(\lambda)}M\subset J_\lambda$. Suppose first
that (a) holds; it follows that the natural map
$$
M/I^nM\to\lim_{\lambda\in\Lambda}M/(I^nM+J_\lambda)
$$
is injective for every $n\in\N$.
Taking into account example \ref{ex_lim_interchange}(ii),
we deduce that both of the induced maps
$$
M\to\lim_{n\in\N}M/I^nM
\qquad\text{and}\qquad
\lim_{n\in\N}M/I^nM\to
\lim_{(\lambda,n)\in\Lambda\times\N}M/(I^nM+J_\lambda)
$$
are injective. However, notice that the subset
$\Sigma\subset\Lambda\times\N$ of all $(n,\lambda)$
such that $n\geq n(\lambda)$, is cofinal, hence
the composition of these two maps is an isomorphism.
We conclude that both these maps are bijective, and
the assertion follows.

Next, suppose that (b) holds, and let $a_1,\dots,a_r$
be a finite system of generators for $I$. Clearly $M$
is separated for the $I$-adic topology. Now, suppose
$x_\bullet:=(x_n~|~n\in\N)$ is a sequence of elements
of $M$ such that $x_n-x_m\in I^nM$ for every $n,m\in\N$
with $m\geq n$. The sequence $x_\bullet$ converges for
the topology $\cT_M$ to an element $x\in M$, and it
remains to check that $x$ is also the limit of $x_\bullet$
for the $I$-adic topology of $M$. To this aim, for every
$n\in\N$, let
$$
S_n:=\{(j_1,\dots,j_r)\in\N^r~|~j_1+\cdots+j_r=n\}
$$
and for every $\sigma:=(j_1,\dots,j_r)\in\N^r$ set
$a^\sigma:=a_1^{j_1}\cdots a_r^{j_r}$.
For every $k,l\in\N$ we may find inductively a system of
elements $(y_{l,\sigma}~|~\sigma\in S_k)$ of $M$ such that
\begin{itemize}
\item
$x_{k+l}-x_k=\sum_{\sigma\in S_k}y_{l,\sigma}a^\sigma$
\item
$y_{l',\sigma}-y_{l,\sigma}\in I^lM$ whenever $l'\geq l$ and
for every $\sigma\in\N^r$.
\end{itemize}
Thus, the sequence $(y_{l,\sigma}~|~l\in\N)$ converges
in the topology $\cT_M$ to some element $y_\sigma\in M$
for every $\sigma\in\N^r$, and we get
$x-x_k=\sum_{\sigma\in S_k}y_\sigma a^\sigma$; so this
difference lies in $I^kM$ for every $k\in\N$, whence
the contention.
\end{proof}

\begin{proposition}\label{prop_f-adics}
Let $(A,\cT)$ be any f-adic ring. We have :
\begin{enumerate}
\item
If $A_0$ is any ring of definition of $A$, then the
topology of $A_0$ is adic and admits a finitely
generated ideal of adic definition.
\item
A subring of $A$ is a ring of definition of $A$ if
and only if it is open and bounded in $A$.
\item
The following conditions are equivalent :
\begin{enumerate}
\item
The topology $\cT$ is linear.
\item
The topology $\cT$ is adic and admits a finitely generated
ideal of adic definition.
\end{enumerate}
\end{enumerate}
\end{proposition}
\begin{proof} Fix an open subring $A_1$ and a finitely
generated ideal $I$ of $A_1$, such that $\cT$ induces the
$I$-adic topology on $A_1$. 

(i): Since $A_0$ is open in $A$, it follows that there exists
an integer $n>0$ with $I^n\subset A_0$; set $J:=I^nA_0$.
If $K$ is any open ideal of $A_0$, we may likewise
find $m\in\N$ such that $I^{mn}\subset K$, so that
$J^m\subset K$ and $J^m$ is open in $A_0$. The assertion
follows easily.

(ii): Clearly, every ring of definition of $A$ is open
and bounded. Conversely, suppose that $A_0$ is open and
bounded. As in the foregoing, we find an integer $n>0$ such
that $I^n\subset A_0$ and we set $J:=I^nA$. Let now $U$ be
any open neighborhood of zero in $A_0$; since $A_1$ admits
a fundamental system of open neighborhoods of zero consisisting
of additive subgroups, we may assume that $U$ is an additive
subgroup of $A_0$. Since $A_0$ is bounded, there exists
$m\in\N$ such that $x\cdot a\in U$ for every $x\in I^m$
and every $a\in A_0$; then clearly $J^mA_0\subset U$.
Lastly, it is easily seen that $J^k$ is open in $A_0$
for every $k\in\N$, whence (ii).

(iii): If $\cT$ is linear, then $A$ is bounded in itself,
so the assertion follows directly from (ii).
\end{proof}

\begin{corollary}\label{cor_f-adics}
Let $(A,\cT)$ be any f-adic ring. We have :
\begin{enumerate}
\item
If $A_0$ and $A_1$ are rings of definition of $A$, then
the same holds for $A_0\cap A_1$ and $A_0\cdot A_1$.
\item
Let $B$ (resp. $C$) be a bounded (resp. open) subring
of $A$, with $B\subset C$. Then there exists a ring
of definition $A_0$ of $A$ with $B\subset A_0\subset C$.
\item
$(A,\cT)$ is adic if and only if $A$ is bounded in
the topology $\cT$.
\item
$A^\circ$ is the filtered union of all the subrings of
definition of $A$.
\end{enumerate}
\end{corollary}
\begin{proof}(i): The assertion for $A_0\cap A_1$ is
immediate from proposition \ref{prop_f-adics}(ii).
The assertion for $A_0\cdot A_1$ follows likewise,
taking into account that $A$ has a fundamental system
of open neighborhoods of zero consisting of additive
subgroups : details left to the reader.

(ii): Let $A_1$ be any ring of definition of $A$;
by proposition \ref{prop_f-adics}(ii), the subring
$C_1:=A_1\cap C$ is a ring of definition of $A$,
and we may take $A_0:=B\cdot C_1$.

(iii): If $(A,\cT)$ is adic, then clearly $A$ is
bounded. Conversely, if $A$ is bounded, (ii) implies
that $A$ is a subring of definition of $(A,\cT)$, and
then $(A,\cT)$ is adic, by proposition \ref{prop_f-adics}(i).

(iv): By remark \ref{rem_someth-on-bdd-in-Z-lin}(iii),
we know that $A^\circ$ is the filtered union of the
family $\cF$ of bounded subrings of $A$; notice that
any subring containing an open subring of $A$ is also
open in $A$. Since the family of bounded and open subring
of $A$ is not empty (proposition \ref{prop_f-adics}(ii)),
we conclude that this family is cofinal in $\cF$, whence
the contention.
\end{proof}

\begin{corollary}\label{cor_Tate}
Let $(A,\cT_A)$ be any topological ring, $B\subset A$ an
open subring, and endow $B$ with the topology $\cT_B$
induced from $\cT_A$. The following holds :
\begin{enumerate}
\item
$(A,\cT_A)$ is f-adic if and only if the same holds for
$(B,\cT_B)$.
\item
Suppose that $(A,\cT_A)$ is a Tate ring, and $B$ is a
ring of definition of $A$. Then we have :
\begin{enumerate}
\item
$B$ contains an element of $A^\times\cap A^{\circ\circ}$.
\item
Let $s\in B\cap A^\times\cap A^{\circ\circ}$ be any element.
Then $sB$ is an ideal of adic definition for $B$, and
$B_s=A$.
\end{enumerate}
\item
Conversely, let $C$ be any ring, $f\in C$ any element,
$D$ the image of\/ $C$ in the localization $C_f$, and
endow $C_f$ with the unique $D$-linear topology $\cT_f$
such that $\{f^nD~|~n\in\N\}$ is a fundamental system of
neighborhoods of\/ $0$ in $C_f$. Then $(C_f,\cT_f)$ is a
Tate ring.
\end{enumerate}
\end{corollary}
\begin{proof}(i): Directly from the definition we see that
if $(B,\cT_B)$ is f-adic, the same holds for $(A,\cT_A)$.
Conversely, if $(A,\cT_A)$ is f-adic, let $A_0\subset A$
be any subring of definition; by corollary \ref{cor_f-adics}(ii)
we may find a subring of definition $A_1$ of $A$ contained
in $A_0\cap B$; especially, $A_1$ is an open subring of $B$
whose topology is adic with a finitely generated ideal of
adic definition (proposition \ref{prop_f-adics}(i)), whence
the assertion.

(ii.a) is clear. Next, for every $a\in A$ there exists
$n\in\N$ such that $s^na\in B$, so $B_s=A$. Also, the
map $a\mapsto s^na$ is an automorphism of $(A,\cT_A)$,
so $s^nB$ is open for every $n\in\N$. Lastly, since
$B$ is bounded in $A$, for every open neighborhood $U$
of $0\in A$ there exists $n\in\N$ such that $s^nB\subset U$,
whence (ii.b).

(iii): Explicitly, the open subsets of $(C_f,\cT_f)$
are the arbitrary unions of subsets of the form $a+f^mD$,
where $a\in C_f$ is arbitrary and $m\in\N$ is any integer;
indeed, it is easily seen that the intersection of any
two subsets of this type is either empty, or equal to
one of them, so such arbitrary unions define a topology
$\cT_f$, and it is easily seen that $(C_f,\cT_f)$ is a
topological ring (details left to the reader). A simple
inspection then shows that $(C_f,\cT_f)$ is a Tate ring. 
\end{proof}

\begin{lemma}\label{lem_5.3.8}
Let $A$ be any c-adic, complete and separated topological
ring, $I,J\subset A$ two ideals, with $I$ finitely generated.
We have :
\begin{enumerate}
\item
$I$ is open in $A$ if and only if the same holds for $I^c$.
\item
If $I$ is an ideal of c-adic definition for $A$, the
following holds :
\begin{enumerate}
\item
$A$ is adic and $I$ is an ideal of adic definition for $A$.
\item
$J$ is open in $A$ if and only if the same holds for $J^c$.
\end{enumerate}
\end{enumerate}
\end{lemma}
\begin{proof} Assertion (i) is a special case of
\cite[Lemma 5.3.5(ii) and 5.3.8(i)]{Ga-Ra}, and
(ii.a) follows directly from (i). Lastly, suppose
that $J^c$ is open in $A$; by (ii.a) there exists
$n\in\N$ such that $I^n\subset J^c\subset J+I^{n+1}$.
Letting $M:=(I^n+J)/J$, it follows that $IM=M$, and
since $A$ is complete and separated for its $I$-adic
topology, $I$ lies in the Jacobson radical of $A$
(remark \ref{rem_someth-on-bdd-in-Z-lin}(v));
therefore $M=0$ by Nakayama's lemma, {\em i.e.}
$I^n\subset J$, so $J$ is open in $A$.
\end{proof}

\begin{lemma}\label{lem_tag-reinstated}
Let $A$ be a complete and separated adic topological
ring, $I\subset A$ an ideal of adic definition, $M$
an $A$-module whose $I$-adic topology $\cT_M$ is
separated, $N\subset M$ a finitely generated submodule,
and $N^c$ the topological closure of $N$ in $(M,\cT_M)$.
We have :
\begin{enumerate}
\item
If the $I$-adic topology $\cT_N$ of $N^c$ agrees
with the one induced by $\cT_M$, then $N=N^c$.
\item
$N$ is open in $(M,\cT_M)$ if and only if the same holds
for $N^c$.
\end{enumerate}
\end{lemma}
\begin{proof}(i): By assumption, $IN^c$ is an open submodule
of $N^c$ for the topology induced by $(M,\cT_M)$, hence
$N^c=N+IN^c$, and the assertion follows from \cite[Th.8.4]{Mat}.

(ii): If $N$ is open in $M$, then $N=N^c$ and the assertion is
clear. Conversely, if $N^c$ is open in $M$, then $I^nM\subset N^c$
for some $n\in\N$, hence $\cT_N$ agrees with the topology
induced from $\cT_M$; then we get again $N=N^c$, by (i).
\end{proof}

\begin{definition}\label{def_c-adic-and-adic}
Let $(A,\cT)$ and $(A',\cT')$ be two topological rings,
and $f:A\to A'$ a ring homomorphism.

(i)\ \
We say that $f$ is {\em c-adic} (resp. {\em adic})
if $\cT$ and $\cT'$ are linear, and the family of ideals
$$
((f(I)\cdot A')^c~|~\text{$I$ open ideal in $A$})
\qquad
\text{(resp.\ \
$(f(I)\cdot A'~|~\text{$I$ open ideal in $A$})$\ )}
$$
is a fundamental system of open neighborhoods of zero
in $A'$.

(ii)\ \
Suppose that $A$ and $A'$ are f-adic. Then we say that
$f$ is {\em f-adic} if there are rings of definition
$A_0$ of $A$ and $A'_0$ of $A'$ such that $f(A_0)\subset A'_0$
and the restriction $f_{|A_0}$ is adic.

(iii)\ \
Let $\bQ$ be any property of ring homomorphisms
({\em e.g.} ``of finite type'', ``flat'', ``\'etale''
and so on). We say that $f$ is an {\em adically $\bQ$ morphism}
(resp. a {\em c-adically $\bQ$ morphism}) if $f$ is
adic (resp. c-adic) and for every open ideal $I\subset A$
the induced map $A/I\to B/(IB)^c$ is a $\bQ$ ring homomorphism.
\end{definition}

\begin{lemma}\label{lem_f-adics}
Let $f:(A,\cT)\to(A',\cT')$ be any ring homomorphism
of topological rings.

{\em(i)}\ \
If $f$ is c-adic, the following holds :
\begin{enumerate}
\alphaenu
\item
$f$ is a continuous ring homomorphism.
\item
If the topology $\cT$ is c-adic, the same holds for $\cT'$,
and if $I$ is any ideal of c-adic definition for $A$, then
$f(I)A'$ is an ideal of c-adic definition for $A'$.
\item
If the topologies $\cT$ and $\cT'$ are f-adic,
linear and complete, then $f$ is adic.
\end{enumerate}

{\em(ii)}\ \ 
If the topologies $\cT$ and $\cT'$ are adic,
the following conditions are equivalent :
\begin{enumerate}
\alphaenu
\item
$f$ is adic.
\item
There exists an ideal of adic definition $I$ of $A$ such
that $f(I)\cdot A'$ is an ideal of adic definition of $A'$.
\item
For every ideal $I$ of adic definition of $A$, the
ideal $f(I)\cdot A'$ is of adic definition for $A'$.
\end{enumerate}

{\em(iii)}\ \
Suppose that $f$ is an f-adic map of f-adic topological
rings. Then the following holds :
\begin{enumerate}
\alphaenu
\item
$f$ is continuous and bounded ({\em i.e.} if $S$ is
any bounded subset of $A$, the subset $f(S)$ is bounded
in $A'$). Especially, $f(A^\circ)\subset A'^\circ$.
\item
For every ring of definition $A_0$ of $A$ and $A'_0$
of $A'$ such that $f(A_0)\subset A'_0$, and every ideal
of adic definition $I$ of $A_0$, the ideal $f(I)\cdot A'_0$
is of adic definition for $A'_0$.
\item
For every ring of definition $A_0$ of $A$ and every
open subring $B\subset A'$ such that $f(A_0)\subset B$,
there exists a ring of definition $A'_0$ of $A'$ such
that $f(A_0)\subset A'_0\subset B$.
\item
If $B\subset A$, $B'\subset A'$ are open subrings with
$f(B)\subset B'$, then the restriction $f_{|B}:B\to B'$
is f-adic, for the topologies of $B$ and $B'$ induced by
$\cT$ and $\cT'$.
\end{enumerate}

{\em(iv)}\ \
If $f$ is a continuous open map, and the topologies
$\cT$ and $\cT'$ are f-adic, then $f$ is f-adic.

{\em(v)}\ \
Suppose that $f(A)$ is open in $A'$, that $\cT$ is a linear
topology and that $\cT'$ is adic. Then $f$ is adic if and
only if it is a continuous open map.
\end{lemma}
\begin{proof}(i.a) is obvious, and (i.c) follows easily
from proposition \ref{prop_f-adics}(iii) and lemma
\ref{lem_5.3.8}(i). To show (i.b), let us first remark,
quite generally :

\begin{claim}\label{cl_image-and-closure}
Let $f:X\to Y$ be any continuous map between arbitrary
topological spaces, and $S,T\subset X$ two subsets such
that $S^c=T^c$. Then $(fS)^c=(fT)^c$.
\end{claim}
\begin{pfclaim} It suffices to show that $fS\subset(fT)^c$.
However, clearly we have $T\subset f^{-1}((fT)^c)$,
therefore $S\subset T^c\subset f^{-1}((fT)^c)$, whence the
contention.
\end{pfclaim}

\begin{claim}\label{cl_exchange-pow-clo}
Let $B$ be any topological ring, $J_1,J_2\subset B$
two additive subgroups. Denote $J_1J_2$ the additive
subgroup generated by
$(j_1j_2~|~j_1\in J_1;\ j_2\in J_2)$, and define likewise
$J_1^cJ_2^c$. Then
$$
(J_1J_2)^c=(J^c_1J^c_2)^c
\qquad
\text{for every $n\in\N$}.
$$
\end{claim}
\begin{pfclaim} Since $J_1J_2\subset J_1^cJ_2^c$, we have
$(J_1J_2)^c\subset(J^c_1J^c_2)^c$. For the converse, it suffices
to show that $J^c_1J^c_2\subset(J_1J_2)^c$. However, for every
$k\in\N$, consider the map
$$
\mu_k:(B^{\oplus 2})^{\oplus k}\to B
\qquad
(b_{ij}~|~i=1,2;\ j=1,\dots,k)\mapsto
\sum_{t=1}^kb_{1t}b_{2t}.
$$
Then $J^c_1J_2^c=
\bigcup_{k\in\N}\mu_k((J^c_1\oplus J_2^c)^{\oplus k})$.
Now, $\mu_k$ is clearly continuous for the product topology
on $(B^{\oplus 2})^{\oplus k}$, and $(J^c_1\oplus J_2^c)^{\oplus k}$
is the topological closure of $(J_1\oplus J_2)^{\oplus k}$ in
this topology, so $\mu_k((J^c_1\oplus J_2^c)^{\oplus k})$ is
contained in the topological closure in $B$ of
$\mu_k((J_1\oplus J_2)^{\oplus k})$ (claim \ref{cl_image-and-closure}),
which in turns is contained in $(J_1J_2)^c$, whence the contention.
\end{pfclaim}

Now, let $I$ be any ideal of c-adic definition of
$A$; by assumption, for every open ideal $I'$ of $A'$
there exists $n\in\N$ such that $(f((I^n)^c)A')^c\subset I$.
Set $J:=f(I)A'$; then $J$ is open in $B$, and by claims
\ref{cl_exchange-pow-clo} and \ref{cl_image-and-closure}
we deduce that
$$
(f((I^n)^c)A')^c=((f((I^n)^c))^cA')^c=((f(I^n))^cA')^c=
(f(I^n)A')^c=(J^n)^c
$$
whence the contention.

(ii): If $I_1,I_2$ are any two ideals of adic definition
of $A$, then there exist $n,m\in\N$ such that
$I_1^n\subset I_2^m\subset I_1$. The assertion is an
easy consequence : details left to the reader.

(iii.a): The continuity of $f$ is obvious. Suppose that
$S\subset A$ is a bounded subset, and let $U'\subset A'$
be any open neighborhood of zero; we have to exhibit
an open neighborhood of zero $V'$ in $A'$ such that
$f(S)\cdot V'\subset U'$. To this aim, let $A_0$ and
$A'_0$ be rings of definition of $A$ and $A'$, such
that $f(A_0)\subset A'_0$, and $I\subset A_0$ an ideal
of adic definition such that $IA'_0$ is an ideal of
adic definition of $A'_0$; we may assume that $U$ is
an open ideal of $A'_0$, and we may find $n,m\in\N$
such that $f(I^n)\subset U'$ and $S\cdot I^m\subset I^n$.
Thus, $f(S)\cdot f(I^m)\subset U'$, so $V':=I^mA'_0$ will do.

(iii.b): Let $A_1$ and $A'_1$ be rings of definitions
of $A$ and $A'$ such that $f(A_1)\subset A'_1$, and
$J\subset A_1$ an ideal of adic definition such that
$K:=f(J)\cdot A'_1$ is an ideal of adic definition of
$A'_1$. Then there exist $n,m\in\N$ such that
$J^n\subset I$ and $K^m\subset A'_0$, so that
$K^{n+m}=f(J^n)\cdot K^m\subset f(I)\cdot A'_0$;
especially, $f(I)\cdot A'_0$ is open in $A'_0$.
Moreover, since $I$ is topologically nilpotent
and finitely generated in $A_0$, it is easily
seen that $f(I)A'_0$ is topologically nilpotent
in $A'_0$; thus $f(I)A'_0$ is an ideal of adic
definition of $A'_0$, as stated.

(iii.c) follows easily from (iii.a) and corollary
\ref{cor_f-adics}(ii).

(iii.d): By corollary \ref{cor_Tate}(i), both $B$ and $B'$
are f-adic rings. Let $A_0\subset A$ and $A'_0\subset A'$
be subrings of definition with $f(A_0)\subset A'_0$; hence
$B_0:=A_0\cap B$ and $B'_0:=A'_0\cap B'$ are subrings of
definition of $B$ and $B'$ respectively, with
$f(B_0)\subset B'_0$ (proposition \ref{prop_f-adics}(ii));
then the assertion follows from (iii.b).

(iv): Let $A_0$ and $A'_0$ be rings of definition of $A$
and $A'$; after replacing $A_0$ by $A_0\cap f^{-1}A'_0$
we may assume that $f(A_0)\subset A'_0$ (proposition
\ref{prop_f-adics}(ii)); by assumption, $f(A_0)$ is
open in $A'$, so we may further replace $A'_0$ by
$f(A'_0)$, and assume from start that $f$ restricts
to an open surjective map $A_0\to A'_0$. In this case,
let $I\subset A_0$ be any ideal of adic definition;
it follows easily that $f(I)$ is an ideal of adic
definition for $A'_0$, so $f$ is adic.

(v)\ \
If $f$ is a continuous open map, then for every pair of
open ideals $J\subset A'$ and $I\subset A$ with $f(I)\subset J$,
we have $f(I)\cdot A'\subset J$, and clearly $f(I)\cdot A'$
is open in $A'$. This shows that $f$ is adic.

Conversely, if $f$ is adic and $f(A)$ is open in $A'$,
let $J\subset f(A')$ be an open ideal of $A'$. For every
open ideal $I\subset A$, the ideal $f(I)\cdot A'$ is open
in $A'$, hence the same holds for $f(I)\cdot J$; but
$f(I)\cdot J\subset f(I)$, so $f(I)$ is open for every
such $I$, and thus $f$ is open.
\end{proof}

\begin{example}\label{ex_f-adic-quotient}
Let $(A,\cT)$ be a topological ring, $f:A\to A'$ a surjective
ring homomorphism, and endow $A'$ with the topology $\cT'$
induced by $A$ via $f$. We have :

(i)\ \
It is easily seen that a subset $U'\subset A'$ is open if and
only if $U'=f(U)$ for some open subset $U$ of $A$. It follows
that $(A',\cT')$ is a topological ring.

(ii)\ \
Moreover, if $(A,\cT)$ is c-adic, the same holds for $(A',\cT')$,
and $f:(A,\cT)\to(A',\cT')$ is c-adic. Indeed, let $I\subset A$
be an ideal of c-adic definition, and set $J:=f(I)$; by (i) we
know that $f((I^n)^c)$ is open in $A$' for every $n\in\N$, and
then claim \ref{cl_image-and-closure} implies that
$f((I^n)^c)=(J^n)^c$, so the latter is open in $A'$ for every
$n\in\N$, again by (i), and the system of such ideals is a
fundamental system of open neighborhoods of $0$ in $A'$.

(iii)\ \
Likewise, if $(A,\cT)$ is adic (resp. f-adic), the same holds
for $(A',\cT')$, and $f$ is adic (resp f-adic).
\end{example}

\begin{definition}\label{def_deja-vu}
For any f-adic ring $A$, set
$$
X_A:=\Spec\,A
\qquad
X_A^{\circ\circ}:=\Spec\,A/A^{\circ\circ}A.
$$
We call $X^{\circ\circ}_A$ the {\em non-analytic locus\/} of $X_A$
and its complement $X_A\setminus X^{\circ\circ}_A$ the
{\em analytic locus}. 
\end{definition}

\begin{lemma}\label{lem_deja-vu}
Let $A$ and $B$ be two f-adic rings, $f:A\to B$
a continuous ring homomorphism, and set $\phi:=\Spec\,f$.
We have :
\begin{enumerate}
\item
$X_A^{\circ\circ}=\{\fp\in X_A~|~\text{$\fp$ is open in $A$}\}$.
\item
$\phi$ restricts to a map $X_B^{\circ\circ}\to X_A^{\circ\circ}$.
\item
If $f$ is an injective and open map, $\phi$ restricts to an
isomorphism of schemes :
$$
X_B\setminus X_B^{\circ\circ}\isom X_A\setminus X_A^{\circ\circ}.
$$
\item
The homomorphism $f$ is f-adic if and only if
$\phi^{-1}(X^{\circ\circ}_A)=X_B^{\circ\circ}$.
\item
Let $J\subset A$ be any ideal. Then $J$ is open if and only
if\/ $\Spec\,A/J\subset X^{\circ\circ}_A$.
\end{enumerate}
\end{lemma}
\begin{proof}(i) and (ii) are clear.

(v): If $J$ is open, then the same holds for every prime
ideal $\fp$ containing $J$; in light of (i) it then follows
that $A^{\circ\circ}\subset\fp$ for every such $\fp$, so
the radical of $J$ also contains $A^{\circ\circ}$, {\em i.e.}
$\Spec\,A/J\subset X^{\circ\circ}_A$. Conversely, suppose
that the radical of $J$ contains $A^{\circ\circ}$, and pick
any subring of definition $A_0\subset A$, and a finitely
generated ideal $I\subset A_0$ of adic definition; it
follows easily that $J$ contains $I^n$ for every sufficiently
large $n\in\N$, so $J$ is open.

Next, let us check that if $f$ is adic, then
$\phi^{-1}(X^{\circ\circ}_A)=X_B^{\circ\circ}$. Indeed, pick
subrings of definition $A_0\subset A$, $B_0\subset B$ such
that $f(A_0)\subset B_0$, and an ideal of adic definition
$I$ of $A_0$; then $J:=f(I)\cdot B_0$ is an ideal of adic
definition of $B_0$ (lemma \ref{lem_f-adics}(iii.b)).
It is easily seen that $X_A^{\circ\circ}=\Spec\,A/IA$ and
$X_B^{\circ\circ}=\Spec\,B/JB$, whence the stated
condition.

(iii): Notice first that $f$ is f-adic in the current
situation, so the foregoing already shows that
$\phi^{-1}(X^{\circ\circ}_A)=X_B^{\circ\circ}$. The
latter means that $X_B\setminus X^{\circ\circ}_B=
\bigcup_{s\in A^{\circ\circ}}\Spec\,B_s$. Now, let
$\fp\in X_A\setminus X_A^{\circ\circ}$ be any
non-open prime ideal of $A$; then there exists
$s\in A^{\circ\circ}$ such that $s\notin\fp$, and it
suffices to check that the localization $f_s:A_s\to B_s$
is an isomorphism. However, $f_s$ is obviously injective;
let $b\in B$ be any element; since $A$ is open in $B$,
there exists $n\in\N$ such that $s^nb\in A$, so $f_s$
is surjective as well.

(iv): By the foregoing, we may assume that
$\phi^{-1}(X^{\circ\circ}_A)=X_B^{\circ\circ}$, and we check
that $f$ is adic. Indeed, pick a subring of
definition $B_0\subset B$; since $f$ is continuous,
$f^{-1}B_0$ is open in $A_0$, so we may find a subring
of definition $A_0$ of $A$ such that $A_0\subset f^{-1}B_0$
(corollary \ref{cor_f-adics}(ii)). Let $f_0:A_0\to B_0$
be the restriction of $f$, and set $\phi_0:=\Spec\,f_0$;
we get a commutative diagram
$$
\xymatrix{ X_B \ar[r]^-\phi \ar[d] & X_A \ar[d] \\
X_{B_0} \ar[r]^-{\phi_0} & X_{A_0}
}$$
whose vertical arrows are induced by the corresponding
inclusion maps. The assumption implies that
$\phi^{-1}(X_A\setminus X^{\circ\circ}_A)=X_B\setminus X^{\circ\circ}_B$;
together with (iii), we deduce that
$\phi_0^{-1}(X_{A_0}\setminus X^{\circ\circ}_{A_0})=
X_{B_0}\setminus X^{\circ\circ}_{B_0}$, which is equivalent
to $\phi_0^{-1}(X^{\circ\circ}_{A_0})=X^{\circ\circ}_{B_0}$. On
the other hand, it suffices to show that $f_0$ is an
adic ring homomorphism, so we may replace $A$ and $B$
by $A_0$ and $B_0$, and assume from start that $A$ and
$B$ are adic topological rings. However, our assumption
on $\phi^{-1}(X^{\circ\circ}_A)$ means that the radical of
$f(A^{\circ\circ})\cdot B$ equals $B^{\circ\circ}$; now, pick
any finitely generated ideal of adic definition $I$ of $A$
(resp. $J$ of $B$); then the radical of $f(I)\cdot B$ also
equals $B^{\circ\circ}$, and according to (v), the latter
condition implies that $f(I)\cdot B$ is open in $B$, so
$f$ is adic.
\end{proof}

\begin{proposition}\label{prop_top-on-opens-fadic-case}
With the notation of definition {\em\ref{def_deja-vu}},
let $U\subset X_A$ be a quasi-compact open subset
containing the analytic locus. We have :
\begin{enumerate}
\item
There exists a unique ring topology $\cT_U$ on
$A_U:=\cO_{X_A}(U)$ such that $(A_U,\cT_U)$ is f-adic
and the restriction map $\rho_U:A\to A_U$ is open.
\item
Suppose moreover that $A$ admits a complete, separated
and noetherian ring of definition. Then $(A_U,\cT_U)$ is
complete and separated.
\item
Let $\phi:A\to B$ be any morphism of f-adic topological
rings, such that the image of\/ $\Spec\,\phi$ lies in $U$.
Then the resulting map $\phi_U:(A_U,\cT_U)\to B$ is continuous.
\end{enumerate}
\end{proposition}
\begin{proof}(i): Clearly, there exists a unique topology
$\cT_U$ on $A_U$ such that $\rho_U$ is open and $(A_U,\cT_U)$
is a topological group for its additive group structure, and
we need to check that $(A_U,\cT_U)$ is a topological ring.
Let $I$ be a finitely generated ideal of adic definition
for a ring of definition $A_0$ of $A$; the assertion comes
down to the following. For every $a\in A_U$ and every
$n\in\N$ there exists $k\in\N$ such that
$$
a\cdot\rho_U(I^k)\subset\rho_U(I^n).
$$
Now, pick a finite system of generators $(f_1,\dots,f_r)$
of $I$, and notice that $A^{\circ\circ}$ is the radical of
the ideal $IA^\circ$ in the subring $A^\circ$ of $A$
(remark \ref{rem_someth-on-bdd-in-Z-lin}(iv)), hence
\set\begin{equation}\label{eq_equal-loci}
X_A^{\circ\circ}=\Spec\,A/IA.
\end{equation}
Let $j:U\to X_A$ be the inclusion map; then $j_*\cO_U$ is a
quasi-coherent $\cO_{X_A}$-module (\cite[Ch.I, Cor.9.2.2]{EGAI}),
and on the other hand, from \eqref{eq_equal-loci} and our
assumption on $U$ we see that $\Spec\,A_{f_i}\subset U$. Hence
$$
A_{f_i}=\cO_{X_A}(\Spec\,A_{f_i})=
j_*\cO_U(\Spec\,A_{f_i})=j_*\cO_U(X_A)_{f_i}=A_U[f_i^{-1}]
\qquad
\text{for $i=1,\dots,r$}.
$$
Thus, there exists $m\in\N$ such that $f_i^ma\in\rho_U(A)$
for every $i=1,\dots,r$. Then there also exists $t\in\N$
such that $\rho_U(I^t)\cdot f^m_ia\subset\rho_U(I^n)$ for
$i=1,\dots,r$, and the assertion follows easily. It is
then also clear that the topology $\cT_U$ is f-adic, and
$\rho_U(A_0)$ is a ring of definition of $A_U$, with $\rho_U(I)$
as ideal of adic definition.

(ii): Let $A_0\subset A$ be a noetherian, complete and separated
subring of definition; it suffices to check that $A_0\cap\Ker\,\rho_U$
is a closed ideal, but this holds, by \cite[Th.8.10(i)]{Mat}.

(iii): Let $H\subset B$ be any open additive subgroup; then
$\phi^{-1}_UH$ contains the open subgroup $\rho_U(\phi^{-1}H)$,
so it is open, whence the contention.
\end{proof}

\begin{lemma}\label{lem_still-c-adic}
Let $f:(A,\cT_A)\to(B,\cT_B)$ be a c-adic morphism
of topological rings. Then:
\begin{enumerate}
\item
The completion map $(A,\cT_A)\to(A^\wedge,\cT^\wedge_A)$
is c-adic, and $A$ is c-adic if and only if the same holds
for $A^\wedge$.
\item
The separated completion
$f^\wedge:(A^\wedge,\cT_A^\wedge)\to(B^\wedge,\cT^\wedge_B)$
of $f$ is c-adic.
\item
Let $\bQ$ be any property of ring homomorphisms. If $f$
is a c-adically $\bQ$ morphism, the same holds for 
$f^\wedge$ (see definition {\em\ref{def_c-adic-and-adic}(iii)}).
\item
If the topology $\cT_A$ is $I$-adic for some finitely
generated ideal $I\subset A$, the topology $\cT^\wedge_A$
is $IA^\wedge$-adic.
\end{enumerate}
\end{lemma}
\begin{proof}(i): Let $(I_\lambda~|~\lambda\in\Lambda)$
be any fundamental system of open ideals of $A$, and
set $I^\wedge_\lambda:=(I_\lambda A^\wedge)^c$ for every
$\lambda\in\Lambda$. From remark
\ref{rem_completion-of-topring}(ii) we see that
$(I^\wedge_\lambda~|~\lambda\in\Lambda)$ is a
fundamental system of open neighborhoods of zero in
$A^\wedge$, so $\cT^\wedge_A$ is a linear topology, and the
completion map $A\to A^\wedge$ is c-adic. In view of lemma
\ref{lem_f-adics}(i.b), it follows that if $A$
is c-adic, the same holds for $A^\wedge$. For the
converse, suppose that $A^\wedge$ is c-adic, and pick
an ideal of c-adic definition of the form $(IA^\wedge)^c$,
where $I$ is an open ideal of $A$; for every $n\in\N$
we let $(I^n)^c$ (resp. $(I^n)^\wedge$) be the topological
closure of $I^n$ in $A$ (resp. in $A^\wedge$); then
$$
(I^n)^c=A\cap(I^n)^\wedge
\qquad
\text{for every $n\in\N$}.
$$
Now, let $J\subset A$ be any open ideal; the topological
closure $J^\wedge$ of $J$ in $A^\wedge$ is open in $A^\wedge$,
so there exists $n\in\N$ such that $(I^n)^\wedge\subset J^\wedge$,
and therefore $(I^n)^c\subset A\cap J^\wedge=J$, which says
that $A$ is c-adic with $I$ as ideal of c-adic definition.

(ii): We know already that $\cT^\wedge_A$ and $\cT^\wedge_B$
are both linear topologies. Next, by assumption, for every
open ideal $J\subset B$ there exists $\lambda\in\Lambda$
such that $(I_\lambda B)^c\subset J$; it follows that
$(I_\lambda^\wedge B^\wedge)^c=(I_\lambda B^\wedge)^c
\subset(JB^\wedge)^c$. Likewise, for every
$\lambda\in\Lambda$ there exists an open ideal
$J$ of $B$ such that $J\subset(I_\lambda B)^c$,
and therefore $(JB^\wedge)^c\subset
(I_\lambda^\wedge B^\wedge)^c$. Since the family
$((JB^\wedge)^c~|~\text{$J$ open in $B$})$ is a
fundamental system of open ideals of $B$, it
follows that the same holds for the family
$((I_\lambda^\wedge B^\wedge)^c~|~\lambda\in\Lambda)$,
whence the contention.

(iii): This follows easily from (ii), taking into
account the natural isomorphisms
$$
A^\wedge/(IA^\wedge)^c\isom A/I
\qquad\text{and}\quad
B^\wedge/(JB^\wedge)^c\isom B/J
$$
for every open ideals $I$ of $A$ and $J$ of $B$, as
well as the equality :
$((IA^\wedge)^cB^\wedge)^c=(IB^\wedge)^c$, which identify
the map $A^\wedge/(IA^\wedge)^c\to B^\wedge/(IB^\wedge)^c$
induced by $f^\wedge$ with the map $A/I\to B/(IB)^c$
induced by $f$.

(iv) follows easily from remark
\ref{rem_completion-of-topring}(ii,iv).
\end{proof}

\begin{proposition}\label{prop_complete-f-adic}
Let $A$ be an f-adic topological ring, $A_0$ an open
subring of $A$, and $A_0^\wedge$, $A^\wedge$ the
separated completions of $A_0$ and respectively $A$.
We have :
\begin{enumerate}
\item
$A^\wedge$ is an f-adic topological ring.
\item
If $A_0$ is a ring of definition of $A$, then $A_0^\wedge$
is a ring of definition for $A^\wedge$.
\item
The natural map $j:A\otimes_{A_0}A_0^\wedge\to A^\wedge$ is a
ring isomorphism.
\end{enumerate}
\end{proposition}
\begin{proof} By corollary \ref{cor_not-in-Bourbaki}(i),
the natural map $A_0^\wedge\to A^\wedge$ is an open immersion,
so (i) and (ii) follow from lemma \ref{lem_still-c-adic}(iv).

(iii): Suppose first that $A_0$ is a ring of definition
of $A$. Since $A_0$ is open in $A$, the topology induced
by $A$ on $B:=A/A_0$ via the projection $A\to B$ is
discrete, so the completion map $B\to B^\wedge$ is
an isomorphism of topological groups, and there follows
a commutative ladder :
$$
\xymatrix{ 0 \ar[r] & A_0^\wedge \ar[r] \ddouble &
A\otimes_{A_0}A_0^\wedge \ar[r] \ar[d]_j &
B\otimes_{A_0}A_0^\wedge \ar[r] \ar[d]^{j_B} & 0 \\
0 \ar[r] & A_0^\wedge \ar[r] & A^\wedge \ar[r] & B \ar[r] & 0
}$$
whose bottom row is short exact, by proposition
\ref{prop_replaces-Mat-Th.8.1}(v). It follows that the top
row is exact as well, and to conclude, it suffices to check
that $j_B$ is an isomorphism. To this aim, let us write $B$
as the colimit of the filtered system
$(B_\lambda~|~\lambda\in\Lambda)$ of its finitely generated
$A_0$-submodules; clearly it suffices to show that the
restriction $j_\lambda:B_\lambda\otimes_{A_0}A_0^\wedge\to B_\lambda$
of $j_B$ is an isomorphism for every such $\lambda$.
However, for every $\lambda\in\Lambda$ there exists
an open ideal $I_\lambda\subset A_0$ such that
$I_\lambda B_\lambda=0$, so $j_\lambda$ factors as the
composition of the isomorphism
$B_\lambda\otimes_{A_0}A^\wedge_0\isom
B_\lambda\otimes_{A_0}A^\wedge_0/I_\lambda A^\wedge_0$ and the
natural identifications
$B_\lambda\otimes_{A_0}A_0^\wedge/I_\lambda A^\wedge_0\isom
B_\lambda\otimes_{A_0}A_0/I_\lambda\isom B_\lambda$, whence
the contention, in this case. Next, if $A_0$ is an
arbitrary open subring of $A$, then $A_0$ is an f-adic
ring (corollary \ref{cor_Tate}(i)); we let $A_1$ be
a ring of definition of $A_0$, and we notice that
$A_1$ is also a ring of definition for $A$ (proposition
\ref{prop_f-adics}(ii)). By the foregoing, both natural
maps $A_0\otimes_{A_1}A_1^\wedge\to A_0^\wedge$ and
$A\otimes_{A_1}A_1^\wedge\to A^\wedge$ are isomorphisms,
so the same follows for $j$.
\end{proof}

\begin{proposition}\label{prop_f-adic-push-out}
Let $A,B,B'$ be three f-adic topological rings and $f:B\to A$,
$g:B\to B'$ two f-adic ring homomorphisms. Then we have :
\begin{enumerate}
\item
There exists a unique f-adic ring topology $\cT_{A'}$ on
$A':=B'\otimes_BA$ such that the resulting diagram (with
$f':=\one_{B'}\otimes_Bf$ and $g_A:=g\otimes_B\one_A$)
\set\begin{equation}\label{eq_both-f-adic}
{\diagram
B \ar[r]^-f \ar[d]_g & A \ar[d]^{g_A} \\
B' \ar[r]^-{f'} & (A',\cT_{A'})
\enddiagram}
\end{equation}
is cocartesian in the category $\Z\tdu\TopAlg_{\Z\tdu\mathrm{lin}}$
(notation of definition {\em\ref{def_top-ring}(v)}).
Moreover, both $f'$ and $g_A$ are f-adic ring homomorphisms.
\item
Suppose additionally that $f$ is open. Then the same holds
for $f'$, and \eqref{eq_both-f-adic} is cocartesian in the
category $\Z\tdu\TopAlg$.
\item
In the situation of {\em (ii)}, suppose moreover that $B$ is a
ring of definition of $A$, and $g$ induces an isomorphism
$g^\wedge:B^\wedge\isom B'^\wedge$ on separated completions; then
\begin{enumerate}
\item
$g_A$ induces an isomorphism $g_A^\wedge:A^\wedge\isom A'^\wedge$
on separated completions
\item
$A'^\circ$ (resp. $A'{}^{\circ\circ}$) is the image of
$B'\otimes_BA^\circ$ (resp. of $B'\otimes_BA^{\circ\circ}$)
in $A'$.
\end{enumerate}
\end{enumerate}
\end{proposition}
\begin{proof}(i): In light of corollary \ref{cor_f-adics}(i)
and lemma \ref{lem_f-adics}(iii.b) we may find subrings of
definitions $B_0\subset B$, $A_0\subset A$ and $B'_0\subset B'$
such that $f$ and $g$ restrict to adic ring homomorphisms
$f_0:B_0\to A_0$ and $g_0:B_0\to A'_0$. Let also $I_0\subset B_0$
be an ideal of adic definition; denote by $\bar A{}'_0\subset A'$
the image of the induced map $A'_0:=B'_0\otimes_{B_0}A'_0\to A'$,
and endow $\bar A{}'_0$ with its $I_0\bar A{}'_0$-adic topology.
There exists a unique group topology $\cT'$ on $A'$ such that
$\bar A{}'_0$ is an open subgroup, and we claim that $\cT'$ is
a ring topology. This comes down to the following assertion.
For every $a'\in A'$ and every $n\in\N$ there exists $m\in\N$
such that
$$
a'\cdot I_0^m\bar A{}'_0\subset I_0^n\bar A{}'_0.
$$
However, we may write $a'$ as a finite sum of terms of
the form $b'\otimes a$, with $a\in A$ and $b'\in B'$,
so we are easily reduced to the case where $a'=b'\otimes a$
for such $a$ and $b'$. But since $f$ and $g$ are f-adic,
we may find, for every $n\in\N$, integers $m_1,m_2\in\N$
such that $a\cdot I_0^{m_1}A_0\subset I^nA_0$, and
$b'\cdot I_0^{m_2}B'_0\subset B'_0$; obviously $m:=m_1+m_2$
will do. It is now clear that $f'$ and $g_A$ are both
f-adic ring homomorphisms for this topology $\cT_{A'}$ on
$A'$. Lastly, let $(C,\cT_C)$ be any topological ring whose
topology is $\Z$-linear, $h:B'\to C$ and $k:A\to C$ two
continuous ring homomorphisms such that $k\circ f=h\circ g$;
there follows a unique ring homomorphism $l:A'\to C$ such
that $l\circ i=h$ and $l\circ g_A=k$. To conclude, we have
to check that $l$ is continuous, and it suffices to show
that the restriction $\bar l_0:\bar A{}'_0\to C$ is
continuous. But the topology of $\bar A{}'_0$ agrees with
the quotient topology induced from the $I_0A'_0$-adic topology
$\cT_{A'_0}$ on $A'_0$, so we reduced to checking that the
induced map $l_0:(A'_0,\cT_{A'_0})\to C$ is continuous in a
neighborhood of $0$. Thus, let $U$ be a neighborhood of $0$
in $C$; since $\cT_C$ is $\Z$-linear, we may assume that $U$ is
a subgroup of $C$. Moreover, there exists an open neighborhood
$U'$ of $0$ in $C$ such that $U'\cdot U'\subset U$. Then we
may find $n\in\N$ such that $k(I^n_0A_0),h(I^n_0B'_0)\subset U'$,
so that $l_0(I^{2n}_0A'_0)\subset U$, whence the contention.

(ii): Indeed, if $f$ is open, we may replace $A_0$ by
$f(B_0)$ in the foregoing, in which case the induced
map $B'_0\to A'_0$ is surjective and open; since the
same holds for the projection $A'_0\to\bar A{}'_0$, we
see that $f'$ is open. Next, let $C$ be any topological
ring, and $h,k$ two continuous ring homomorphisms as in
the foregoing. We have to show that the induced map
$l:A'\to C$ is continuous, and we are again reduced
to showing that the same holds for its restriction
$\bar l_0:\bar A{}'_0\to C$. But we have just seen that
the topology of $\bar A{}'_0$ agrees with the quotient
topology induced by the surjective map $f'_0:B'_0\to\bar A{}'$
(the restriction of $f'$); thus $\bar l_0$ is continuous
if and only if the same holds for $\bar l_0\circ f'_0$,
and the latter is none else than the restriction of $h$.

(iii): Indeed, under these assumptions, we get a
cocartesian diagram of rings
$$
\xymatrix{ B' \ar[rr]^-{f'} \ar[d]_j & & A' \ar[d] \\
B'^\wedge \ar[r]^-{g^{\wedge-1}} & B^\wedge \ar[r] &
B^\wedge\otimes_{B'}A' \ar[r]^-\sim &
B^\wedge\otimes_BA \ar[r]^-\sim & A^\wedge
}$$
where $j$ is the completion map, and the righmost
isomorphism on the bottom row is provided by proposition
\ref{prop_complete-f-adic}(iii). Since the natural map
$B^\wedge\to A^\wedge$ is injective, it follows especially
that $\Ker\,f'\subset\Ker\,j$, and therefore $f'$ induce
an isomorphism $B'{}^\wedge\isom\bar A{}_0^{\prime\wedge}$
on separated completions. Invoking again proposition
\ref{prop_complete-f-adic}(iii) we deduce a commutative
diagram
$$
\xymatrix{
B^\wedge\otimes_BA \ar[r] \ar[d] & B'^\wedge\otimes_{B'}A' \ar[r] &
\bar A{}_0^{\prime\wedge}\otimes_{\bar A{}_0'}A' \ar[d] \\
A^\wedge \ar[rr]^-{g_A^\wedge} & & A'^\wedge
}$$
whose top horizontal and vertical arrows are isomorphisms,
so the same holds for the bottom horizontal arrow. Set
$\bar A:=\Img\,g_A$, and endow $\bar A$ with the quotient
topology induced by $A$ via the surjection $A\to\bar A$.
Recall that the topology of $A$ (resp. of $A'$) agrees with
the topology induced by $A^\wedge$ (resp. by $A'^\wedge$) via
the completion map $A\to A^\wedge$ (resp. $A'\to A'^\wedge$);
since we already know that $g^\wedge_A$ is an isomorphism of
topological rings, it follows easily that the topology of
$\bar A$ agrees with the topology induced by $A'$ via the
inclusion map $\bar A\to A'$ (details left to the reader).
Moreover, recall that $g_A$ is f-adic, so lemma
\ref{lem_f-adics}(iii) implies that the image of
$B'\otimes_BA{}^\circ$ lies in $A'^\circ$; for the converse, say
that $x\in A'^\circ$, and write $x=\sum_{i=1}^nb_i\otimes a_i$
for some $n\in\N$ and elements $a_i\in A$, $b_i\in B'$, for
$i=1,\dots,n$. Let $I\subset B$ be any ideal of definition;
we may find $m\in\N$ such that $a_iI^m\subset A^{\circ\circ}$
for every $i=1,\dots,n$, and since $g^\wedge$ is an isomorphism,
we may also find $b'_i\in B$, $b''_i\in I^mB'$ such that
$b_i=g(b'_i)+b''_i$ for $i=1,\dots,n$.
Clearly $b''_i\otimes a_i$ lies in the image of
$B'\otimes_BA{}^{\circ\circ}$ for every $i=1,\dots,n$,
so $y:=1\otimes(\sum_{i=1}^na_ib'_i)\in\bar A\cap A'^\circ
\subset\bar A{}^\circ$, and it remains only to observe
that $\bar A{}^\circ=g_A(A^\circ)$, because the kernel
of the surjection $A\to\bar A$ is contained in the kernel
of the completion map $A\to A^\wedge$ (details left to the
reader). Likewise, if $x\in A'^{\circ\circ}$, then
$y\in\bar A\cap A'^{\circ\circ}=\bar A{}^{\circ\circ}=g(A^{\circ\circ})$,
which shows that $A'{}^{\circ\circ}$ lies in the image of
$B'\otimes_BA{}^{\circ\circ}$, and the converse inclusion
is obvious.
\end{proof}

\sset\subsubsection{}\label{subsec_polynomial-top-rings}
Let $A$ be any topological ring; for any two subsets
$S,T\subset A$ we denote by $S\cdot T$ the additive
subgroup of $A$ generated by $\{st~|~s\in S,\ t\in T\}$.
We define inductively
$$
T^0:=\{1\}
\qquad\text{and}\qquad
T^{n+1}:=T\cdot T^n
\qquad
\text{for every $n\in\N$}.
$$
Now, let $(T_\lambda~|~\lambda\in\Lambda)$ be any (small)
family of subsets of $A$; then, for every $\nu\in\N^{(\Lambda)}$
the additive subgoup
$$
T^\nu:=\prod_{\lambda\in\Lambda}T^{\nu(\lambda)}_\lambda
$$
is well defined. We consider the polynomial $A$-algebra
$$
A[X_\bullet]:=A[X_\lambda~|~\lambda\in\Lambda]
$$
and we also set $X^\nu:=\prod_{\lambda\in\Lambda}X_\lambda^{\nu(\lambda)}$
for every such $\nu$. Furthermore, let $\cU_A$ be the set
of all open neighborhood $U$ of $0$ in $A$, and set
$$
U[T_\bullet X_\bullet]:=
\Bigl\{\sum_{\nu\in\N^{(\Lambda)}}a_\nu X^\nu\in A[X_\bullet]~|~
\text{$a_\nu\in T^\nu\cdot U$ for every $\nu\in\N^{(\Lambda)}$}\Bigr\}
\qquad
\text{for every $U\in\cU_A$}.
$$

\begin{proposition}\label{prop_polynomial-top-rings}
In the situation of \eqref{subsec_polynomial-top-rings},
suppose that $A$ is an f-adic ring, and $T_\lambda\cdot A$
is an open ideal of $A$ for every $\lambda\in\Lambda$,
and let $A_0\subset A$ be a subring of definition. Then we have :

{\em (i)}\ \
There exists a unique ring topology $\cT_{A,T_\bullet}$ on
$A[X_\bullet]$ for which the family
$$
\cU_{A[X_\bullet]}:=
(U[T_\bullet X_\bullet]~|~U\in\cU_A)
$$
is a fundamental system of open neighborhoods of\/ $0$.

{\em (ii)}\ \
The topological ring
$$
A[X_\bullet]_T:=(A[X_\bullet],\cT_{A,T_\bullet})
$$
is f-adic, and $A_0[T_\bullet X_\bullet]$ is a subring of definition
of $A[X_\bullet]$. Moreover, if $I\subset A_0$ is an ideal of adic
definition, then $I[T_\bullet X_\bullet]$ is an ideal of adic
definition for $A_0[T_\bullet X_\bullet]$.

{\em (iii)}\ \
The structure map $h:A\to A[X_\bullet]$ is f-adic for the
topology $\cT_{A,T_\bullet}$, and the subset
$\{tX_\lambda~|~\lambda\in\Lambda,\ t\in T_\lambda\}$ is
power bounded in $A[X_\bullet]$.

{\em (iv)}\ \
Moreover, for every f-adic ring $B$, every family
$(b_\lambda~|~\lambda\in\Lambda)$ of elements of $B$,
and every continuous ring homomorphism $f:A\to B$
such that
$F:=\{f(t)\cdot b_\lambda~|~\lambda\in\Lambda,\ t\in T_\lambda\}$
is power bounded in $B$, there exists a unique continuous
ring homomorphism $g:(A[X_\bullet],\cT_{A,T_\bullet})\to B$ with
$g\circ h=f$ and $g(X_\lambda)=b_\lambda$ for every
$\lambda\in\Lambda$.
\end{proposition}
\begin{proof}(i): Notice that every element of $\cU_{A[X_\bullet]}$
is an additive subgroup of $A[X_\bullet]$; it follows easily that
there exists a unique group topology on $A[X_\bullet]$ that
admits $\cU_{A[X_\bullet]}$ as a fundamental system of neighborhoods
of $0$. Moreover, for $U\in\cU_A$ we may find $V\in\cU_A$ with
$V\cdot V\subset U$, whence
$V[T_\bullet X_\bullet]\cdot V[T_\bullet X_\bullet]\subset
U[T_\bullet X_\bullet]$. Therefore, it remains only to check
that for every $P\in A[X_\bullet]$ and every $U\in\cU_A$
there exists $V\in\cU_A$ such that
$$
P\cdot V[T_\bullet X_\bullet]\subset U[T_\bullet X_\bullet].
$$
To this aim, notice first that, since $A$ is f-adic and
$T_\lambda\cdot A$ is open in $A$ for every $\lambda\in\Lambda$,
then $T^\nu\cdot U\in\cU_A$ for every $\nu\in\N^{(\Lambda)}$
and every $U\in\cU_A$ (details left to the reader).
Say that $P=\sum_{\nu\in\N^{(\Lambda)}}a_\nu X^\nu$; since
$a_\nu=0$ except for finitely many $\nu\in\N^{(\Lambda)}$,
we may find $V\in\cU_A$ small enough that
$a_\nu\cdot V\in T^\nu\cdot U$ for every $\nu\in\N^{(\Lambda)}$.
It is easily seen that such $V$ will do.

(ii): From the foregoing, it is clear that
$A_0[T_\bullet X_\bullet]$ is an open subring of $A[X_\bullet]$
and $I[T_\bullet X_\bullet]$ is an ideal of $A_0[T_\bullet X_\bullet]$.
Moreover, $I^n[T_\bullet X_\bullet]=I[T_\bullet X_\bullet]^n$ for
every $n\in\N$, so the topology $\cT_{A,T_\bullet}$ induces the
$I[T_\bullet X_\bullet]$-adic topology on $A_0[T_\bullet X_\bullet]$.

(iii): The continuity of $h$ is obvious; moreover, $h$
restricts to a map $A_0\to A_0[T_\bullet X_\bullet]$, and
$I\cdot A_0[T_\bullet X_\bullet]=I[T_\bullet X_\bullet]$, so
$h$ is f-adic.
Lastly, we have $tX_\lambda\in A_0[T_\bullet X_\bullet]$ for
every $\lambda\in\Lambda$ and $t\in T_\lambda$, whence the
assertion.

(iv): Obviously there exists a unique map of $A$-algebras
$g:A[X_\bullet]\to B$ such that $g(X_\lambda)=b_\lambda$ for
every $\lambda\in\Lambda$, and it suffices to check that
$g$ is continuous. To this aim, for any open additive
subgroup $U$ of $B$ pick $V\in\cU_B$ such that
$F^n\cdot V\subset U$ for every $n\in\N$. Set also
$W:=f^{-1}V$; then $g(W[T_\bullet X_\bullet])\subset U$,
whence the contention.
\end{proof}

\sset\subsubsection{Modules of finite type over topological rings}
\label{subsec_fg-mods-top-rings}
Let $(A,\cT_A)$ be a topological ring, and $M$ any $A$-module
of finite type. We choose any surjective $A$-linear map
$\pi:A^{\oplus n}\to M$ (for some $n\in\N$) and we endow
$A^{\oplus n}:=A\times\cdots\times A$ with the product of the
topologies of its direct factors, and $M$ with the topology
$\cT^A_M$ induced by $A^{\oplus n}$ via the map $\pi$.

\begin{lemma}\label{lem_cantop-on-fg-mods}
With the notation of \eqref{subsec_fg-mods-top-rings}, the
following holds :
\begin{enumerate}
\item
The topology $\cT^A_M$ is independent of the choice of $\pi$.
\item
$\pi:(A,\cT_A)^n\to(M,\cT^A_M)$ is an open map.
\end{enumerate}
\end{lemma}
\begin{proof}(i): Let $\pi':A^{\oplus m}\to M$ be any other
$A$-linear surjection; we define
$\pi'':A^{\oplus n+m}=A^{\oplus n}\oplus A^{\oplus m}\to M$ as the
unique $A$-linear map whose restriction to $A^{\oplus n}\oplus 0$
agrees with $\pi$ and whose restriction to $0\oplus A^{\oplus m}$
agrees with $\pi'$. It suffices to check that the topologies on
$M$ induced by $A^{\oplus n}$ and $A^{\oplus m}$ via $\pi$ and
respectively $\pi'$ agree with the topology $\cT''$ induced by
$A^{\oplus n+m}$ via $\pi''$.
Thus, we are reduced to comparing $\cT^A_M$ and $\cT''$. Now,
let $e_1,\dots,e_{n+m}$ be the standard basis of $A^{\oplus n+m}$,
and for each $i=1,\dots,m$ choose $f_i\in A^{\oplus n}\oplus 0$
such that $\pi''(f_i)=\pi''(e_{i+n})$; we consider the $A$-linear
automorphism
$$
\omega:A^{\oplus n+m}\isom A^{\oplus n+m}
\qquad
e_j\mapsto\left\{\begin{array}{ll}
                    e_j & \text{for $j=1,\dots,n$} \\
                    e_j-f_{j-n} & \text{for $j=n+1,\dots,n+m$}.
                 \end{array}\right.
$$
Clearly $\omega:(A,\cT_A)^{n+m}\isom(A,\cT_A)^{n+m}$ is a
homeomorphism, and the restriction of $\pi''\circ\omega$
to $A^{\oplus n}\oplus 0$ still agrees with $\pi$. We may then
replace $\pi''$ by $\pi''\circ\omega$ and assume from
start that $0\oplus A^{\oplus m}\subset\Ker\,\pi''$. Now, let
$U\subset M$ be any subset; we have $U\in\cT^A_M$ if and
only if $\pi^{-1}U$ is open in $(A,\cT_A)^n$, and the
latter holds if and only if $(\pi^{-1}U)\times A^{\oplus m}$
is open in $(A,\cT_A)^{n+m}$, which in turn is equivalent
to $U\in\cT''$, whence the contention.

(ii): Let $U\subset A^{\oplus n}$ be any open subset and set
$U':=U+\Ker\,\pi:=\{u+v~|~u\in U,\ v\in\Ker\,\pi\}$; then
$U'$ is open in $A^{\oplus n}$ and $\pi^{-1}(\pi U)=U'$, whence
the contention.
\end{proof}

\begin{definition} In the situation of
\eqref{subsec_fg-mods-top-rings}, we call $\cT^A_M$ the
{\em canonical topology} of the $A$-module $M$ of finite
type.
\end{definition}

\begin{proposition}\label{prop_cantop-on-fg-mods}
Let $(A,\cT_A)$ be a topological ring, $M_1,\dots,M_k$ a finite
sequence of $A$-modules of finite type, $(N,\cT_N)$ a topological
$A$-module and $\beta:M_1\times\cdots\times M_k\to N$ any
$A$-multilinear map. Then we have :
\begin{enumerate}
\item
The canonical topology of the $A$-module
$M_1\times\cdots\times M_k$ is the product of the canonical
topologies $\cT^A_{M_1},\dots,\cT^A_{M_k}$.
\item
$\beta$ is a continuous map
$(M_1,\cT^A_{M_1})\times\cdots\times(M_k,\cT^A_{M_K})\to(N,\cT_N)$.
\end{enumerate}
\end{proposition}
\begin{proof}(i): By a simple induction on $k$, we are easily
reduced to the case where $k=2$. Let then $\cT'$ be the product
of the topologies $\cT^A_{M_1}$ and $\cT_{M_2}$; pick $A$-linear
surjections $\pi_i:A^{\oplus n_i}\to M_i$ for $i=1,2$. It is easily
seen that if $U\in\cT'$, then $(\pi_1\times\pi_2)^{-1}U$ is open in
$(A,\cT_A)^{n_1+n_2}$, hence $U\in\cT^A_{M_1\times M_2}$ (details left
to the reader). Conversely, if $U\in\cT^A_{M_1\times M_2}$, then
$(\pi_1\times\pi_2)^{-1}U$ is a union of subsets of the form
$V_1\times V_2$, where $V_i\subset A^{\oplus n_1}$ is an arbitrary
open subset, for $i=1,2$; therefore
$(\pi_1\times\pi_2)(V_1\times V_2)=(\pi_1V_1)\times(\pi_2V_2)\in\cT''$
since $\pi_1$ and $\pi_2$ are open maps (lemma
\ref{lem_cantop-on-fg-mods}(ii)), and so $U\in\cT''$.

(ii): Pick $A$-linear surjections $\pi_i:A^{\oplus n_i}\to M_i$ for
$i=1,\dots,k$, and set $n:=\sum_{i=1}^kn_i$ and
$\beta':=\beta\circ(\pi_1\times\cdots\times\pi_k)$.
It is easily seen that $\beta':(A,\cT_A)^n\to(N,\cT_N)$
is a continuous map (details left to the reader).
Now, let $U\subset N$ be any open subset; we have
$\beta^{-1}(U)=(\pi_1\times\cdots\times\pi_k)(\beta'^{-1}U)$,
which is open in $M_1\times\cdots\times M_k$, since
$\pi_1\times\cdots\times\pi_k$ is an open map (lemma
\ref{lem_cantop-on-fg-mods}(ii)).
\end{proof}

\begin{corollary}\label{cor_cantop-on-fg-mods}
Let $(A,\cT_A)$ be any topological ring, $B\to A$ a ring
homomorphism, $I\subset A$ an ideal, $M$ and $N$ two
$A$-modules of finite type. The following holds :
\begin{enumerate}
\item
$(M,\cT^A_M)$ is a topological $A$-module.
\item
Every $A$-linear surjection $M\to N$ is a continuous open
map $(M,\cT^A_M)\to(N,\cT^A_N)$.
\item
If the topology $\cT_A$ is $B$-linear (resp. $I$-adic), the
same holds for $\cT^A_M$.
\end{enumerate}
\end{corollary}
\begin{proof} Pick an $A$-linear surjection $\pi:A^{\oplus n}\to M$.

(i): Let $\sigma_M:M\oplus M\to M$ and
$\mu_M:A\times M\to M$ be respectively the addition law and
the scalar multiplication law of $M$; we need to check that
these maps are continuous for the product topologies
$(M,\cT^A_M)\times(M,\cT^A_M)$ and $(A,\cT_A)\times(M,\cT^A_M)$.
However, we have the commutative diagrams
$$
\xymatrix{ A^{\oplus n}\oplus A^{\oplus n}
\ar[r]^-{\sigma_{A^n}} \ar[d]_{\pi\oplus\pi} &
A^{\oplus n} \ar[d]^\pi &
A\times A^{\oplus n} \ar[r]^-{\mu_{A^n}} \ar[d]_{\one_A\times\pi} &
A^{\oplus n} \ar[d]^\pi \\
M\oplus M \ar[r]^-{\sigma_M} & M & A\times M \ar[r]^-{\mu_M} & M
}$$
where $\sigma_{A^n}$ and $\mu_{A^n}$ are the addition law and the
scalar multiplication law of $A^{\oplus n}$, which are clearly
continuous for the product topologies. Now, if $U\in\cT^A_M$, then
$V:=\pi^{-1}U$ is open in $(A,\cT_A)^n$ and $W:=\sigma^{-1}_{A^n}(V)$
is open in $A^{\oplus n}\oplus A^{\oplus n}$; hence
$\sigma^{-1}_M(U)=(\pi\oplus\pi)(W)$ is open in $M\oplus M$, since
$\pi\oplus\pi$ is open (lemma \ref{lem_cantop-on-fg-mods}(ii)).
Likewise one checks the continuity of $\mu_M$.

(ii): The continuity of $f:M\to N$ follows from (i) and proposition
\ref{prop_cantop-on-fg-mods}(ii). Next, we easily reduce to checking
that the map $f\circ\pi:(A,\cT_A)^n\to(N,\cT^A_N)$ is open; but this
map is surjective, so the assertion follows from lemma
\ref{lem_cantop-on-fg-mods}(i).

(iii): Suppose that $(J_\lambda~|~\lambda\in\Lambda)$ is
a fundamental system of open neighborhoods of $0\in A$
consisting of $B$-submodules; since $\pi$ is open, the
family $(\pi(J^{\oplus n}_\lambda)~|~\lambda\in\Lambda)$ is
a fundamental system of open neighborhoods of $0\in M$
(details left to the reader); in view of (i), the assertion
follows. Likewise, if $\cT_A$ is the $I$-adic topology, then
$(I^kA^{\oplus n}~|~k\in\N)$ is a fundamental system of open
neighborhoods of $0\in A^{\oplus n}$, so $(I^kM~|~k\in\N)$ is
a fundamental system of open neighborhoods of $0\in M$,
{\em i.e.} $\cT^A_M$ is the $I$-adic topology of $M$.
\end{proof}

\begin{remark}\label{rem_cantops-on-fin-algs}
(i)\ \
Let $(A,\cT_A)$ be a topological ring and $B$ any finite
$A$-algebra. In light of proposition
\ref{prop_cantop-on-fg-mods}(ii) and corollary
\ref{cor_cantop-on-fg-mods}(i), it is easily seen that
$(B,\cT^A_B)$ is a topological $A$-algebra; moreover,
if $\cT_A$ is the $I$-adic topology (for some ideal
$I\subset A$), then $\cT^A_B$ is the $IB$-adic topology
(corollary \ref{cor_cantop-on-fg-mods}(ii)).

(ii)\ \
Let $N$ be any finite $B$-module. On the one hand, we may
endow $N$ with its canonical topology $\cT^B_N$ as $B$-module
of finite type (that is, relative to the topological
ring $(B,\cT^A_B)$); on the other hand, by restriction of
scalars, $N$ is also an $A$-module of finite type, so we
have as well the canonical topology $\cT^A_N$. We claim
that these topologies coincide. Indeed, pick an $A$-linear
surjection $\pi:A^{\oplus n}\to B$ and a $B$-linear surjection
$\pi':B^{\oplus m}\to N$; there follows an $A$-linear surjection
$\pi^{\oplus m}:A^{\oplus nm}\to B^{\oplus m}$ and by virtue of
proposition \ref{prop_cantop-on-fg-mods}(i) the product
topology of $B^{\oplus m}=(B,\cT^A_B)\times\cdots\times(B,\cT^A_B)$
agrees with the topology induced by $A^{\oplus nm}$ via
$\pi^{\oplus m}$. Then the topology $\cT^B_N$ is also induced
by $\pi'\circ\pi^{\oplus m}$, whence the contention.

(iii)\ \
In the situation of (i), suppose that the topology $\cT_A$
is $f$-adic; let $A_0\subset A$ be a subring of definition,
$I_0\subset A_0$ a finitely generated ideal of adic definition
and $(x_\lambda~|~\lambda\in\Lambda)$ a finite system of
generators of $I_0$. Say that $B=\sum^n_{i=1}Ab_i$ for some
elements $b_1,\dots,b_n\in B$ with $b_1:=1$, and let
$\pi:A^{\oplus n}\to B$ be the $A$-linear surjection such
that $\pi(a_1,\dots,a_n)=\sum_{i=1}^na_ib_i$ for every
$a_1,\dots,a_n\in A$. We have a system
$(c_{ijk}~|~i,j,k=1,\dots,n)$ of elements of $A$ such that
$$
b_ib_j=\sum_{k=1}^nc_{ijk}b_k
\qquad
\text{for every $i,j=1,\dots,n$}.
$$
After replacing $I_0$ by $I_0^k$ for $k\in\N$ large enough,
we may assume that $I_0c_{ijk}\in A_0$ for every
$i,j,k=1,\dots,n$. Let $B_0\subset B$ be the $A_0$-submodule
generated by
$\{1\}\cup\{b_ix_\lambda~|~\lambda\in\Lambda,\, i=2,\dots,n\}$.
It is easily seen that $B_0$ is a subring of $B$, and moreover
$J:=\pi(I_0^{\oplus n})\subset B_0$, so $B_0$ is open in $B$
(lemma \ref{lem_cantop-on-fg-mods}(ii)). Furthermore, a simple
calculation shows that $J$ is an ideal of $B_0$, and clearly
$I_0B_0\subset J$. Lastly, the family
$(\pi((I_0^{k+1})^{\oplus n})~|~k\in\N)$ is a fundamental system
of neighborhoods of $0\in B$, and we have
$\pi((I_0^{k+1})^{\oplus n})=I_0^kJ$ for every $k\in\N$, so the
topology of $B_0$ induced by the inclusion into $(B,\cT^A_B)$
is $I_0B_0$-adic. We conclude that the topology $\cT^A_B$ is
f-adic, and the structure map $A\to B$ is f-adic and restricts
to an adic finite ring homomorphism $A_0\to B_0$ between
subrings of definition.
\end{remark}

\subsection{Topologically local and topologically henselian rings}
Localization and henselization are basic techniques in the
study of commutative algebra; in this section we introduce
suitable variants of these constructions for f-adic topological
rings.

\begin{definition}
(i)\ \
Let $\underline A:=A$ be any f-adic ring.
\begin{enumerate}
\item
We say that $A$ is {\em topologically local} if
$A^{\circ\circ}$ lies in the Jacobson radical of $A^\circ$.
\item
We say that $A$ is {\em topologically henselian} if
$(A^\circ,A^{\circ\circ})$ is a henselian pair.
\end{enumerate}
\end{definition}

\begin{proposition}\label{prop_quasi-affinoid}
Let $A$ be any f-adic ring, and $B$ any open subring
of $A$. We have :
\begin{enumerate}
\item
$A$ is topologically local if and only if the same
holds for $B$.
\item
If $A$ is topologically henselian, then $A$ is also
topologically local.
\item
$A$ is topologically henselian if and only if the
same holds for $B$.
\item
Suppose that $A$ is topologically henselian, let $C$
be any finite $A$-algebra, and endow $C$ with its canonical
$A$-module topology $\cT^A_C$ (see remark
{\em\ref{rem_cantops-on-fin-algs}(i)}). Then $(C,\cT^A_C)$
is a topologically henselian f-adic ring.
\end{enumerate}
\end{proposition}
\begin{proof}(i): Indeed, suppose first that $B$ is
topologically local, and let $a\in A^{\circ\circ}$ be
any element; we may find $n\in\N$ such that
$a^{n+1}\in B^{\circ\circ}$, therefore $1-a^{n+1}\in B^{\circ\times}$
and
$$
(1-a)^{-1}=(1+a+\cdots+a^n)/(1-a^{n+1})\in A^\circ
$$
so $a$ lies in the Jacobson radical of $A^\circ$.
Conversely, suppose that $A$ is topologically local,
and let $b\in B^{\circ\circ}$ be any element; then
$1-b\in A^{\circ\times}$, and we may find $n\in\N$
such that $b^{n+1}/(1-b)\in B^\circ$. Therefore
$$
(1-b)^{-1}=b^{n+1}/(1-b)+1+b+\cdots+b^n\in B^\circ
$$
so that $b$ lies in the Jacobson radical of $B^\circ$.

(ii): Quite generally, if $(R,I)$ is a henselian pair,
then $I$ lies in the Jacobson radical of $R$ (details
left to the reader).

(iii): Indeed, suppose first that $B$ is topologically
henselian; then $B$ is topologically local, by (ii), and
therefore the same holds for $A$, by (i); in view of
\cite[Rem.5.1.10(ii)]{Ga-Ra} we then see that $A$ is
topologically henselian if and only if the same holds
for the open subring $C:=B[A^{\circ\circ}]$ (which is f-adic,
by corollary \ref{cor_Tate}(i)).
However, it is easily seen that the inclusion map
$B^\circ\to C^\circ$ is integral, and moreover the ideal
$A^{\circ\circ}$ of $C^\circ$ is the radical of
$B^{\circ\circ}\cdot C^\circ$ (details left to the reader).
Then the assertion follows from \cite[Rem.5.1.10(i,v)]{Ga-Ra}.

In order to prove the converse, consider -- for every
scheme $X$ -- the set $\mathrm{oc}(X)$ of all open and
closed subsets of $X$; clearly any morphism of schemes
$\phi:Y\to X$ induces a mapping
$$
\mathrm{oc}(\phi):\mathrm{oc}(X)\to\mathrm{oc}(Y)
\qquad
Z\mapsto\phi^{-1}Z.
$$
With this notation we have, quite generally :

\begin{claim}\label{cl_open-closed-equalize}
Let $\phi:Y\to X$ be any closed (resp. open) and
surjective morphism of schemes, and $p_1,p_2:Y\times_XY\to Y$
the induced projections. Then the diagram of sets
$$
\xymatrix{
\mathrm{oc}(X) \ar[rr]^-{\mathrm{oc}(\phi)} & &
\mathrm{oc}(Y) \ar@<.5ex>[rr]^-{\mathrm{oc}(p_1)}
\ar@<-.5ex>[rr]_-{\mathrm{oc}(p_2)} & &
\mathrm{oc}(Y\times_XY)
}$$
identifies $\mathrm{oc}(X)$ with the equalizer of the
mappings $\mathrm{oc}(p_1)$ and $\mathrm{oc}(p_2)$.
\end{claim}
\begin{pfclaim} Since $\phi$ is surjective,
$\mathrm{oc}(\phi)$ is injective, and obviously its
image lies in the equalizer of $\mathrm{oc}(p_1)$ and
$\mathrm{oc}(p_2)$. Conversely, let $Z$ be an open and
closed subset of $Y$ such that $p_1^{-1}Z=p_2^{-1}Z$.
Since $\phi$ is surjective, we have
$Z\subset\phi^{-1}\phi(Z)$, and we claim that in fact
$Z=\phi^{-1}\phi(Z)$. Indeed, let $y\in\phi^{-1}\phi(Z)$
be any point; then there exists $z\in Z$ such that
$\phi(y)=\phi(z)$, and we may find $w\in Y\times_XY$
such that $p_1(w)=z$ and $p_2(w)=y$. The condition on
$Z$ then implies that $y\in Z$ as well, as required.
Next, set $Z':=Y\setminus Z$, and notice that $Z'$
lies as well in the equalizer of $\mathrm{oc}(p_1)$
and $\mathrm{oc}(p_2)$, so we have also $Z'=\phi^{-1}\phi(Z')$;
it follows that $\phi(Z)\cap\phi(Z')=\emptyset$, and since
$\phi$ is a closed (resp. open) mapping, we conclude that
$\phi(Z)$ is an open and closed subset of $X$ whose
preimage in $Y$ equals $Z$, whence the claim.
\end{pfclaim}

Now, if $A$ is topologically henselian, the foregoing
shows that the same holds for $C$, and taking into account
lemma \ref{lem_deja-vu}(iii), to conclude the proof of (iii)
it suffices to remark :

\begin{claim}\label{cl_Jacobson-or-Hensel}
Let $f:R\to S$ be an integral ring homomorphism, $I\subset R$
any ideal. Suppose that the image of the induced morphism of
schemes $\phi:\Spec\,S\to\Spec\,R$ contains
$\Spec\,R\setminus\Spec\,R/I$. Then we have :
\begin{enumerate}
\item
$I$ lies in the Jacobson radical of $R$ if and only
if $IS$ lies in the Jacobson radical of $S$.
\item
The pair $(R,I)$ is henselian if and only if the same
holds for the pair $(S,IS)$.
\end{enumerate}
\end{claim}
\begin{pfclaim} Set $S':=S\times(R/I)$ and let
$\pi:R\to R/I$ be the projection; the pair $(f,\pi)$
determines a unique integral ring homomorphism
$f':R\to S'$, and it is easily seen that the
pair $(S,IS)$ is henselian if and only if the same
holds for the pair $(S',IS')$. Likewise, obviously
$IS$ lies in the Jacobson radical of $S$ if and
only if $IS'$ lies in the Jacobson radical of $S'$.
Moreover, $\Spec\,f'$ is surjective; thus, we may
replace $S$ by $S'$, and assume from start that
$\phi$ is surjective.

(i): Suppose first that $I$ lies in the Jacobson
radical of $R$, and let $\fm\subset S$ be any
maximal ideal; then $\fm':=\phi(\fm)$ is a maximal
ideal of $R$ (\cite[Th.9.4(i)]{Mat}), therefore
$I\subset\fm'$, and hence $IS\subset\fm$, so $IS$
lies in the Jacobson radical of $S$. Conversely,
suppose that the latter condition holds, and let
$\fm'$ be any maximal ideal of $R$; by assumption,
we may find a prime ideal $\fm$ of $S$ with
$\phi(\fm)=\fm'$, and then $\fm$ must be a maximal
ideal, so $IS\subset\fm$, and finally $I\subset\fm'$.

(ii): If $(R,I)$ is henselian, the same holds for $(S,IS)$
(\cite[Rem.5.1.10(v)]{Ga-Ra}), so we may assume that
$(S,IS)$ is henselian, and we show that the same holds
for $(R,I)$. Indeed, let $R\to R'$ be any finite ring
homomorphism; we have to check that the induced
map $R'\to R'/IR'$ restricts to a bijection
$\mathrm{oc}(\Spec\,R')\isom\mathrm{oc}(\Spec\,R'/IR')$.
However, set $S':=R'\otimes_RS$, and notice that the
pair $(S',IS')$ is henselian (again, by
\cite[Rem.5.1.10(v)]{Ga-Ra}), and moreover the induced
map $\Spec\,S'\to\Spec\,R'$ is surjective; thus, we may
replace $R$ by $R'$, and we reduce to showing that the
projection $R\to R/I$ induces a bijection
$\mathrm{oc}(\Spec\,R)\to\mathrm{oc}(\Spec\,R/I)$.
However, $\phi$ is a universally closed morphism
(\cite[Ch.II, Prop.6.1.10]{EGAII}), so claim
\ref{cl_open-closed-equalize} further reduces to
checking that the induced projections $S\to S_0:=S/IS$
and $S\otimes_RS\to S_0\otimes_RS_0$ induce bijections
$$
\mathrm{oc}(\Spec\,S)\to\mathrm{oc}(\Spec\,S_0)
\qquad
\mathrm{oc}(\Spec\,S\otimes_RS)\to
\mathrm{oc}(\Spec\,S_0\otimes_RS_0).
$$
In turns, this is clear, since $(S,IS)$ and
$(S\otimes_RS,IS\otimes_RS)$ are both henselian pairs
(again, by \cite[Rem.5.1.10(v)]{Ga-Ra}).
\end{pfclaim}

(iv): Let $A_0\subset A$ be a subring of definition,
and $I_0\subset A_0$ an ideal of adic definition; by
remark \ref{rem_cantops-on-fin-algs}(iii), we know
that the topology $\cT^A_C$ is f-adic, and $C$ admits
a subring $C_0$ of definition such that the structure
map $A\to C$ restricts to an adic finite ring homomorphism
$A_0\to C_0$. Then the pair $(C_0,I_0C_0)$ is henselian
((\cite[Rem.5.1.10(v)]{Ga-Ra}), and therefore $(C,\cT^A_C)$
is topologically henselian, by (iii).
\end{proof}

\begin{remark}\label{rem_quasi-affinoid}
If $A$ is a topologically local f-adic ring, we have a
homeomorphism
\set\begin{equation}\label{eq_cont-inverse}
1+A^{\circ\circ}\isom 1+A^{\circ\circ}
\quad : \quad
x\mapsto x^{-1}.
\end{equation}
Indeed, say that $a\in A^{\circ\circ}$; then there exists
$n\in\N$ such that
\set\begin{equation}\label{eq_pow-series}
1/(1-a)-(1+a+\cdots+a^n)=a^{n+1}/(1-a)\in A^{\circ\circ}
\end{equation}
and since $A^{\circ\circ}$ is an ideal of $A^\circ$, we deduce
that $1/(1-a)\in 1+A^{\circ\circ}$. It remains only to check
that \eqref{eq_cont-inverse} is continuous on $A^{\circ\circ}$.
However, if $a,b\in A^{\circ\circ}$, we may write
$$
(1-a)^{-1}-(1-a-b)^{-1}=(1-a)^{-1}\cdot(1-1/(1-(1-a)^{-1}b))
$$
so we are reduced to checking the continuity of
\eqref{eq_cont-inverse} at the point $x=1$, and the
latter follows easily from \eqref{eq_pow-series}.
Notice that the same argument proves more precisely
that
$$
(1+U)^{-1}=1+U
$$
for every open additive subgroup $U\subset A^{\circ\circ}$
such that $U\cdot U\subset U$ : details left to the reader.
\end{remark}

\sset\subsubsection{}\label{subsec_localize-f-adic}
Let $A$ be any f-adic ring, $B$ a ring of definition of $A$,
and $I$ a finitely generated ideal of adic definition of $B$.
We let
$$
B_\loc:=(1+I)^{-1}B
\qquad\text{and}\qquad
A_\loc:=A\otimes_BB_\loc
$$
and endow $B_\loc$ with the unique topology such that the
localization map $B\to B_\loc$ is adic, and $A_\loc$ with
the unique f-adic ring topology $\cT_\loc$ such that the
inclusion map $B_\loc\to A_\loc$ is open and the induced
map $A\to A_\loc$ is f-adic (proposition
\ref{prop_f-adic-push-out}(ii)). We claim that $A_\loc$
is topologically local. Indeed, suppose that
$a\in A_\loc^{\circ\circ}$; we may find $n\in\N$ such that
$a^n\in IB_\loc$, so $a^n=b/(1+t)$ for some $b,t\in I$. Therefore,
$1-a^n=(1+t-b)/(1+t)\in B_\loc^\times$, and arguing as in the proof
of proposition \ref{prop_quasi-affinoid}(i) we deduce that
$1-a\in A_\loc^\times$, whence the claim. It follows especially
that $A_\loc$ is independent of the choice of $I$ : indeed,
suppose that $J\subset I$ is another ideal of adic definition
for $B$, and let $A':=(1+J)^{-1}A$; then there exist a
unique map of $A$-algebras $\rho':A'\to A_\loc$, and a natural
isomorphism $A_\loc\isom(1+I)^{-1}A'$ of $A$-algebras that
identifies $\rho'$ with the localization map; however, the
foregoing shows that the image of $1+I$ lies in $A'^\times$,
so $\rho'$ is bijective, and a simple inspection then shows
that $\rho''$ is also an isomorphism of topological rings.
Lastly, let $f:A\to C$ be any continuous ring homomorphism
of f-adic rings, with $C$ topologically local, and let $C_0$
be a ring of definition of $C$, and $J$ an ideal of adic
definition for $C_0$; by the foregoing, we may replace $I$
by a smaller open ideal, and assume that $f(I)\subset J$.
Moreover $C_0$ is topologically local, by proposition
\ref{prop_quasi-affinoid}(i), whence
$f(1+I)\subset 1+J\subset C_0^\times$, so that $f$ factors
uniquely through the localization map $A\to A_\loc$ and a ring
homomorphism $g:A_\loc\to C$. Furthermore, $g$ is a continuous
map : indeed, $g(IB_\loc)=g((1+I)^{-1}I)\subset(1+J)^{-1}J\subset J$,
whence the assertion. Thus, the localization map
$\rho:A\to A_\loc$ is initial in the category of
topologically local topological $A$-algebras with f-adic
topologies; especially, the pair $(A_\loc,\rho)$ is determined
up to unique isomorphism, and we call it the {\em topological
localization} of the f-adic ring $A$.

\sset\subsubsection{}\label{subsec_top-henselization}
For any topological ring $R$, denote by
$$
R\tdu\TopAlg_\mathrm{f-adic}
\qquad\text{and}\qquad
R\tdu\TopHens_\mathrm{f-adic}
$$
the full subcategories of $R\tdu\TopAlg$ whose objects are
respectively the f-adic topological $R$-algebras, and the
topologically henselian f-adic topological $R$-algebras
(definition \ref{def_top-ring}(iii)).
Now, let again $A$, $B$ and $I$ as in \eqref{subsec_localize-f-adic},
denote by $(B^\he_I,I^\he)$ the henselization of the pair $(B,I)$,
endow $B^\he_I$ with its $I^\he$-adic topology, and set
$$
A^\he:=A\otimes_BB^\he_I.
$$
We endow $A^\he$ with the unique topology $\cT^\he_A$ such that
$(A^\he,\cT^\he_A)$ is an f-adic topological ring, the natural map
$B^\he_I\to A^\he$ is open and the natural map $A\to A^\he$ is
f-adic (proposition \ref{prop_f-adic-push-out}(ii)). We may then
state :

\begin{theorem}\label{th_hensel-f-adic}
With the notation of \eqref{subsec_top-henselization},
the following holds :
\begin{enumerate}
\item
$(A^\he,\cT_A^\he)$ is a topologically henselian f-adic ring,
independent of the choice of $B$ and $I$, up to unique
isomorphism of topological $A$-algebras.
\item
The rule : $A\mapsto(A^\he,\cT^\he_A)$ extends to a functor
$$
\Z\tdu\TopAlg_\mathrm{f-adic}\to\Z\tdu\TopHens_\mathrm{f-adic}
$$
that is left adjoint to the forgetful functor (here $\Z$
is endowed with its discrete topology).
\end{enumerate}
\end{theorem}
\begin{proof} Let us consider as well the full subcategories
of $\Z\tdu\TopAlg$ denoted
$$
\Z\tdu\TopAlg_\mathrm{adic}
\qquad\text{and}\qquad
\Z\tdu\TopHens_\mathrm{adic}
$$
whose objects are respectively the adic topological rings,
and the adic topological rings $R$ such that the pair
$(R,R^{\circ\circ})$ is henselian (notice that $R^{\circ\circ}$
is an ideal of $R$).
For every adic ring $R$, we let $(R^\he,(R^{\circ\circ})^\he)$
the henselization of the pair $(R,R^{\circ\circ})$, and we
endow $R^\he$ with the unique ring topology $\cT^\he_R$ such
that the natural map $R\to R^\he$ is adic. We notice :

\begin{claim}\label{cl_first-for-adic}
The rule $R\mapsto(R^\he,\cT^\he_R)$ extends to a functor
$$
\Z\tdu\TopAlg_\mathrm{adic}\to\Z\tdu\TopHens_\mathrm{adic}
$$
that is left adjoint to the forgetful functor.
\end{claim}
\begin{pfclaim} More generally, let $f:R\to S$ and $g:R\to R'$
be two ring homomorphisms, and $I\subset R$ and $J\subset S$
two ideals such that :
\begin{itemize}
\item
$f(I)\subset J$
\item
$J$ lies in the Jacobson radical ideal of $S$
\item
$g$ is unramified, and $g\otimes_RR/I$ is an isomorphism.
\end{itemize}
Then we claim that there exists at most one ring homomorphism
$h:R'\to S$ such that $h\circ g=f$. Indeed, let $h,h':R'\to S$
be any two such maps; according to \cite[Ch.IV, Prop.17.4.6]{EGA4},
the maximal subscheme $U\subset\Spec\,S$ such that
$(\Spec\,h)_{|U}=(\Spec\,h')_{|U}$ is open and closed in
$\Spec\,S$, and on the other hand, $U$ contains $\Spec\,S/J$,
since $g\otimes_RR/I$ is an isomorphism. But $\Spec\,S/J$
meets every open and closed subset of $\Spec\,S$, since
$J$ lies in the Jacobson radical of $S$, whence the assertion.
Next, suppose that the pair $(S,J)$ is henselian, and $g$ is
\'etale; in this case, the map $g\otimes_RS:S\to R'\otimes_RS$
admits a section $s:R'\otimes_RS\to S$
(\cite[Ch.IV, Prop.18.5.4]{EGA4}), and we set
$h:=s\circ(R'\otimes_Rf):R'\to S$. By the foregoing,
$h$ is the unique ring homomorphism such that $h\circ g=f$.
Now, let $(R^\he,I^\he)$ be the henselization of the pair
$(R,I)$, and recall that $R^\he$ is the colimit of a filtered
system $(g_\lambda:R\to R'_\lambda~|~\lambda\in\Lambda)$ of
\'etale ring homomorphisms such that $g_\lambda\otimes_RR/I$
is an isomorphism for every $\lambda\in\Lambda$. Summing
up, we conclude that $f$ factors uniquely through a morphism
of $R$-algebras $R^\he\to S$. Lastly, if $R$ and $S$ are adic
topological rings, then every continuous ring homomorphism
$g:R\to S$ maps $I:=R^{\circ\circ}$ into $J:=S^{\circ\circ}$, and
if $(S,S^{\circ\circ})$ is henselian, the foregoing yields a
unique ring homomorphism $f^\he:R^\he\to S^\he$ extending $f$;
to conclude, it suffices to check that $f^\he$ is continuous.
However, if $J'$ is any ideal of adic definition for $S$, we
may find an ideal $I'$ of adic definition for $R$ such that
$f(I')\subset J'$; by construction, $I'R^\he$ is an ideal
of adic definition of $R^\he$, and $f^\he(I'R^\he)\subset J'$,
whence the contention.
\end{pfclaim}

(i): Let us check first that $A^\he$ is independent of
the choice of $B$ and $I$. Indeed, since $B^{\circ\circ}$
is the radical of $I$, there exists a unique isomorphism
$\phi:B^\he_I\isom B^\he$ of $B$-algebras, where $B^\he$ is as
in the foregoing; it is then clear that the $I^\he$-adic
topology on $B^\he$ agrees with the topology $\cT^\he_B$,
under the isomorphism $\phi$. Let us endow $A\otimes_BB^\he$
with the unique ring topology such that the natural map
$B^\he\to A\otimes_BB^\he$ is open (proposition
\ref{prop_f-adic-push-out}(ii)); it follows that
$\phi\otimes_BA:A^\he\isom A\otimes_BB^\he$ is
an isomorphism of topological $A$-algebras, which shows that
$A^\he$ does not depend on $I$. By the same token, we also
deduce that $A^\he$ is topologically henselian (proposition
\ref{prop_quasi-affinoid}(iii)). Next, let $C\subset A$ be
another ring of definition, $I_C\subset C$ an ideal of definition
of $C$, and $(C^\he,I_C^\he)$ the henselization of the pair
$(C,I_C)$; endow $C^\he$ with its $I^\he_C$-adic topology,
set $A':=A\otimes_CC^\he$, and endow again $A'$ with the
unique ring topology such that the natural map $C^\he\to A'$
is open and the induced map $A\to A'$ is f-adic. We need to
show :

\begin{claim} There exists a unique isomorphism of topological
$A$-algebras $A^\he\isom A'$.
\end{claim}
\begin{pfclaim} In light of corollary \ref{cor_f-adics}(i) we
are easily reduced to the case where $B\subset C$, and since
we have already established that the construction of $C^\he$
and $B^\he$ is independent of the ideals of definition, we may
moreover assume that $I_C=I$ (details left to the reader). In
this situation, set $D:=B^\he\otimes_BC$; the induced ring
homomorphism $B^\he\to D$ identifies $I^\he=B^\he\otimes_BI$ with
the ideal $ID$ of $D$, and therefore the pair $(D,ID)$ is
henselian (\cite[Rem.5.1.10(ii)]{Ga-Ra}). Let us endow
$D$ with its $ID$-adic topology; then $D^{\circ\circ}$ is
the radical of $ID$, and therefore $D$ is topologically
henselian (\cite[Rem.5.1.10(i)]{Ga-Ra}); moreover, the
natural map $C\to D$ is continuous, so by claim
\ref{cl_first-for-adic} it factors uniquely through a
continuous map $C^\he\to D$ of $C$-algebras. We claim that
the latter is an isomorphism of topological rings. Indeed,
let $E$ be any other object of $\Z\tdu\TopHens_\mathrm{adic}$;
any continuous ring homomorphism $B\to E$ admits a
unique continuous extension $B^\he\to E$, so the set
of continuous ring homomorphisms $D\to E$ is in natural
bijection with the set of continuous ring homomorphisms
$C\to E$, and the latter are as well in natural bijection
with the set of continuous ring homomorphisms $C^\he\to E$,
whence the claim (details left to the reader). There
follows an isomorphism of $A$-algebras
$$
\psi:A^\he\isom A\otimes_CD\isom A'
$$
and since the structure maps of these $A$-algebras are
f-adic, it is clear that $\psi$ is an isomorphism of
topological rings. Furthermore, suppose that
$\psi':A^\he\isom A'$ is another isomorphism of topological
$A$-algebras, and set $R:=\psi(B^\he)\cdot\psi'(B^\he)$;
then $R$ is a ring of definition of $A'$ (corollary
\ref{cor_f-adics}(i)), so it is an object of
$\Z\tdu\TopHens_\mathrm{adic}$ (proposition
\ref{prop_quasi-affinoid}(iii)), and the restrictions
$B^\he\to R$ of $\psi$ and $\psi'$ agree on $B$. By
claim \ref{cl_first-for-adic}, it follows that $\psi$
and $\psi'$ agree on $B^\he$, and therefore they coincide,
whence the contention.
\end{pfclaim}

Lastly, let $f:A\to R$ be any continuous ring homomorphism
from $A$ to a topologically henselian f-adic topological ring,
and pick any ring of definition $R_0$ of $R$; after replacing
$B$ by $B\cap f^{-1}R_0$ we may assume that $f$ restricts to
a map $f_0:B\to R_0$ (proposition \ref{prop_f-adics}(ii)),
and since $R_0$ is topologically henselian
(proposition \ref{prop_quasi-affinoid}(iii)), the map
$f_0$ factors uniquely through a continuous ring homomorphism
$f_0^\he:B^\he\to E$ (claim \ref{cl_first-for-adic}). The
datum of $f$ and $f_0^\he$ determines a unique continuous
map of $A$-algebras $A^\he\to E$ (proposition
\ref{prop_f-adic-push-out}(ii)), so the proof is concluded.
\end{proof}

\begin{definition}\label{def_henselize-f-adic}
Let $A$ be any f-adic topological ring. The topological
$A$-algebra $(A^\he,\cT^\he_A)$ provided by theorem
\ref{th_hensel-f-adic} is called the {\em topological
henselization} of $A$.
\end{definition}

\begin{corollary}\label{cor_justify}
Let $A$ be any f-adic topological ring, and denote by
$A_\loc$, $A^\he$, $A^\wedge$ respectively the
topologically localization, the topological henselization,
and the separated completion of $A$. The following holds :
\begin{enumerate}
\item
The localization map $A\to A_\loc$ and the
henselization map $A\to A^\he$ induce isomorphisms
of topological rings on separated completions :
$$
A^{\he\wedge}\xleftarrow{\sim}A^\wedge\isom A_\loc^\wedge.
$$
\item
The inclusion maps $A^\circ\to A$ and $A^{\circ\circ}\to A$
induce natural identifications :
$$
\begin{aligned}
(A^\circ)_\loc\isom\,& (A_\loc)^\circ
& \qquad & &
(A^\circ)^\he\isom\,& (A^\he)^\circ
& \qquad & &
(A^\circ)^\wedge\isom\,& (A^\wedge)^\circ \\
A^{\circ\circ}\otimes_{A^\circ}A^\circ_\loc
\isom\,& (A_\loc)^{\circ\circ}
& \qquad & &
A^{\circ\circ}\otimes_{A^\circ}A^{\circ\he}\isom\,& (A^\he)^{\circ\circ}
& \qquad & &
(A^{\circ\circ})^\wedge\isom\,& (A^\wedge)^{\circ\circ}.
\end{aligned}
$$
\item
Let $R$ be any open subring of $A^\circ$, and denote by
$(R^\he,R^{\circ\circ\he})$ the henselization
of the pair $(R,R^{\circ\circ})$. Then we have :
\begin{enumerate}
\item
There exists a unique ring topology $\cT^\he_R$ on $R^\he$
such that the natural map $R\to R^\he$ is f-adic and
$(R^\he,\cT^\he_R)$ is isomorphic to the topological
henselization of $R$.
\item
Endow $A\otimes_RR^\he$ with the unique f-adic topology such
that the natural map $R^\he\to A\otimes_RR^\he$ is open
(see proposition {\em\ref{prop_f-adic-push-out}(ii)}).
There exists a unique isomorphism of topological $A$-algebras
$A^\he\isom A\otimes_RR^\he$.
\end{enumerate}
\end{enumerate}
\end{corollary}
\begin{proof}(i): Let $B$ be any ring of definition of $A$;
a simple inspection of the construction in
\eqref{subsec_localize-f-adic} shows that the natural map
$$
A':=A\otimes_BB_\loc\to A_\loc
$$
is an isomorphism of topological rings, for the f-adic topology
on $A'$ described in proposition \ref{prop_f-adic-push-out}(ii).
Moreover, the localization map $B\to B_\loc$ clearly
induces an isomorphism $B^\wedge\isom(B_\loc)^\wedge$.
Then the assertion for $A_\loc^\wedge$ follows from
proposition \ref{prop_f-adic-push-out}(iii.a).
Likewise, the natural map $B\to B^\he$ induces an
isomorphism on separated completions, so the same argument
yields the assertion for $A^{\he\wedge}$.

(ii): The assertions for $(A_\loc)^\circ$ and
$(A_\loc)^{\circ\circ}$ follow by the same token,
from proposition \ref{prop_f-adic-push-out}(iii.b).
Likewise, we get the assertions for $(A^\he)^\circ$ and
$(A^\he)^{\circ\circ}$.
In the case of the completion functor, we know that
the natural map $A\otimes_BB^\wedge\to A^\wedge$ is an
isomorphism of topological rings (proposition
\ref{prop_complete-f-adic}(iii)), so we may still
apply proposition \ref{prop_f-adic-push-out}(iii) to deduce
that $(A^\wedge)^\circ$ (resp. $(A^\wedge)^{\circ\circ}$) is the
topological closure in $A^\wedge$ of the image of $A^\circ$
(resp. of $A^{\circ\circ}$).

(iii): Suppose first that $R=A$; especially, $A^\circ=A$.
Then, according to corollary \ref{cor_f-adics}(iv),
the ring $A$ is the filtered union of the system
$(B_\lambda~|~\lambda\in\Lambda)$ of its subrings
of definitions containing $B$, whence an induced
isomorphism of $A$-algebras
$$
A^\he\isom\colim_{\lambda\in\Lambda}B_\lambda\otimes_BB^\he
$$
where $B^\he$ denotes the topological henselization of
$B$. On the other hand, $B$ is also a subring of definition
of $B_\lambda$, and $B^{\circ\circ}_\lambda=B_\lambda\cap B^{\circ\circ}$
for every $\lambda\in\Lambda$. Let
$(B^\he_\lambda,J_\lambda)$ be the henselization of the pair
$(B_\lambda,B^{\circ\circ}_\lambda)$; by theorem
\ref{th_hensel-f-adic}(i) we deduce a unique isomorphism
$$
B_\lambda^\he\isom B_\lambda\otimes_BB^\he
\qquad
\text{for every $\lambda\in\Lambda$}
$$
of topological $B_\lambda$-algebras. Lastly, we have a
unique isomorphism of $A$-algebras
$$
\colim_{\lambda\in\Lambda}B_\lambda^\he\isom R^\he
$$
whence the sought isomorphism $\phi:A^\he\isom R^\he$ of
$A$-algebras. The uniqueness of $\phi$ follows from the
universal property of the henselization functor. If we
may endow $R^\he$ with the topology induced from $A^\he$
via $\phi$, clearly both (iii.a) and (iii.b) hold in this
case.

Next, if $R\subset A^\circ$ is an arbitrary open subring,
choose a ring of definition $B$ of $A$ contained in $R$;
by the foregoing case, $R^\he$ is naturally isomorphic to
the topological henselization of $R$, so that we have unique
isomorphisms
$$
R^\he\isom R\otimes_{A_0}A_0^\he
\qquad
A^\he\isom A\otimes_{A_0}A_0^\he
$$
respectively of topological $R$-algebras and topological
$A$-algebras, such that both natural maps $A^\he_0\to R^\he$
and $A_0^\he\to A^\he$ are open. There follows an isomorphism
of topological $A$-algebras $\phi:A^\he\isom A\otimes_RR^\he$
as sought. The uniqueness of $\phi$ follows again from the
universal property of the topological henselization $A^\he$ of $A$.
\end{proof}

\subsection{Graded structures on topological rings}
Let $\Gamma$ be a monoid, $A:=\bigoplus_{\gamma\in\Gamma}A_\gamma$
a $\Gamma$-graded ring, and $\cT$ a topology on $A$. The
$\Gamma$-graded structure of $A$ is not usually inherited
by the completion $(A,\cT)^\wedge$ of $(A,\cT)$. This observation
motivates the following definition, that introduces a
more flexible notion of graded structure on topological
rings, preserved under completions and related operations.

\begin{definition}\label{def_graded-top}
Let $(A,\cT)$ be a separated topological ring, $\Gamma$
a monoid.

(i)\ \
A {\em $\Gamma$-pre-graded structure} on $(A,\cT)$ is a
datum $(A,\underline B,\Gamma)$ consisting of a subring
$B\subset A$ with a $\Gamma$-graded $\Z$-algebra structure
$\underline B$ on $B$ (definition \ref{def_Gamma-graded-algs}(i))
such that the following holds :
\begin{enumerate}
\alphaenu
\item
The topology induced on $B$ by $\cT$ agrees
with the linear topology defined by a cofiltered system
of graded ideals of $B$ (see definition
\ref{def_Gamma-graded-algs}(iii)).
\item
$B$ is a dense subset of $A$.
\end{enumerate}

(ii)\ \
A {\em $\Gamma$-graded structure} on $(A,\cT)$ is a
$\Gamma$-pre-graded structure $(A,\underline B)$ such
that $\gr_\gamma B$ is a closed subset of $A$, for every
$\gamma\in\Gamma$.

(iii)\ \
Let $(A,\underline B,\Gamma)$ and $(A',\underline B',\Gamma')$
be two topological rings with pre-graded structures.
A {\em morphism of topological rings with pre-graded
structures} $(A,\underline B)\to(A',\underline B')$ is
a pair $(f,\phi)$ consisting of a continuous ring
homomorphism $f:A\to A'$ and a morphism of monoids
$\phi:\Gamma\to\Gamma'$ such that $f(B)\subset B'$,
and such that the restriction $f_{|B}:B\to B'$ induces
a morphism of $\Gamma$-graded $\Z$-algebras (notation
of definition \ref{def_Gamma-graded-algs}(iv))
$$
\underline B\to\Gamma\times_{\Gamma'}\underline B'.
$$
Clearly, we get therefore a category of topological
rings with pre-graded structures :
$$
\pregrTopAlg.
$$
We also have the full subcategory of $\pregrTopAlg$ denoted
$$
\grTopAlg
$$
whose objects are the topological rings with graded structures.
\end{definition}

\begin{remark}\label{rem_graded-top-algs}
Let $(\Gamma,+,0)$ be a monoid, $(A,\underline B)$ a
$\Gamma$-pre-graded structure on the topological ring
$(A,\cT)$. Pick a fundamental system
$(J_\lambda~|~\lambda\in\Lambda)$ of graded open ideals
of $B$, and endow $\gr_\gamma B$ with the topology
$\cT_\gamma$ induced by $\cT$, for every $\gamma\in\Gamma$.

(i)\ \
Since $J_\lambda=\bigoplus_{\gamma\in\Gamma}\gr_\gamma J_\lambda$
for every $\lambda\in\Lambda$, the topology $\cT_\gamma$
agrees with the linear topology defined by the cofiltered
system of $\gr_0B$-submodules
$(\gr_\gamma J_\lambda~|~\lambda\in\Lambda)$, for every
$\gamma\in\Gamma$. Moreover, a graded ideal $I$ of $B$
is closed in the topology of $B$ if and only if
$\gr_\gamma I$ is closed in $\gr_\gamma B$ for every
$\gamma\in\Gamma$. Indeed, since $I\cap\gr_\gamma B=\gr_\gamma I$
for every such $\gamma$, the condition is clearly necessary.
Conversely, notice that the natural projection
$\pi_\gamma:B\to\gr_\gamma B$ is continuous for every
$\gamma\in\Gamma$, hence the same holds for the composition
$\rho_\gamma:B\to\gr_\gamma B/I$ of $\pi_\gamma$
with the quotient map $\gr_\gamma B\to\gr_\gamma B/I$;
now, if $\gr_\gamma I$ is closed in $\gr_\gamma B$, the
topology of $\gr_\gamma B/I$ is separated, hence
$\Ker\,\rho_\gamma$ is closed in $B$, and to conclude it
suffices to remark that $I$ is the intersection of such
kernels, for $\gamma$ ranging over all elements of $\Gamma$.

(ii)\ \
For every $\gamma\in\Gamma$, denote by
$(\gr_\gamma B^\wedge,\cT^\wedge_\gamma)$ the completion
of $(\gr_\gamma B,\cT_\gamma)$. Also, let
$\gr_\gamma J^\wedge_\lambda$ be the topological closure
of $\gr_\gamma J_\lambda$ in $\gr_\gamma B^\wedge$, for
every $\lambda\in\Lambda$ and $\gamma\in\Gamma$, and set
$$
J'_\lambda:=\prod_{\gamma\in\Gamma}\gr_\gamma J^\wedge_\lambda
\qquad
B':=\prod_{\gamma\in\Gamma}\gr_\gamma B^\wedge
\qquad
J''_\lambda:=\bigoplus_{\gamma\in\Gamma}\gr_\gamma J^\wedge_\lambda
\qquad
B'':=\bigoplus_{\gamma\in\Gamma}\gr_\gamma B^\wedge.
$$
Notice that $B''$ (resp. $B'$) is a ring (resp. a
$\gr_0B$-module), and $J''_\lambda$ (resp. $J'_\lambda$)
is an ideal of $B''$ (resp. a $\gr_0B$-submodule of $B'$),
for every $\lambda\in\Lambda$. We endow $B'$ with the
linear topology defined by the cofiltered system of submodules
$(J'_\lambda~|~\lambda\in\Lambda)$. Notice that $B'$ is
complete and separated, and the induced map
$$
j_B:B\to B'
$$
is continuous; more precisely, the topology of $B$ agrees
with the topology induced by $B'$ via $j_B$. Moreover, for
every $\gamma\in\Gamma$, let
$\pi'_\gamma:B'\to\gr_\gamma B^\wedge$ be the projection; since
$J'_\lambda\subset\pi'^{-1}_\gamma(\gr_\gamma J_\lambda)$
for every $\lambda\in\Lambda$, we see that $\pi'_\gamma$
is continuous for the topology $\cT_\gamma$, for every
$\gamma\in\Gamma$.

(iii)\ \
Next, since $B$ is dense in $A$, the inclusion map
$B\to A$ extends uniquely to an isomorphism between
the completion $(A^\wedge,\cT^\wedge)$ of $(A,\cT)$ and
the completion of $B$ (theorem
\ref{th_complete-top-grps}(iii)); it follows that
$j_B$ factors uniquely through a continuous map of
topological $\gr_0B$-modules
$$
j_A:A^\wedge\to B'
\qquad
a\mapsto(a_\gamma~|~\gamma\in\Gamma)
$$
and $j_A$ is injective, since the topology of $B$ is induced
by that of $B'$ (proposition \ref{prop_replaces-Mat-Th.8.1}(i)).
For every $\gamma\in\Gamma$ we call
$\pi^\wedge_\gamma:=\pi'_\gamma\circ j_A:A^\wedge\to\gr_\gamma B^\wedge$ :
$a\mapsto a_\gamma$ the {\em canonical $\gamma$-projection} and set
$$
a_S:=\sum_{\gamma\in S}a_\gamma
\qquad
\text{for every $a\in A^\wedge$ and every finite
subset $S\subset\Gamma$}.
$$

(iv)\ \
For every finite subset $S\subset\Gamma$, let
$i_S:\gr_SB:=\bigoplus_{\gamma\in S}\gr_\gamma B\to A$
be the inclusion map; since $J_\lambda$ is a graded
ideal for every $\lambda\in\Lambda$, it is easily
seen that the topology on $\gr_SB$ induced by $A$
via $i_S$ is the product topology
$\prod_{\gamma\in S}\cT_\gamma$. Hence, $i_S$
extends uniquely to a continuous map
$$
\prod_{\gamma\in S}(\gr_\gamma B^\wedge,\cT^\wedge_\gamma)
=(\gr_SB)^\wedge\to(A^\wedge,\cT^\wedge)
$$
which is still injective (proposition
\ref{prop_replaces-Mat-Th.8.1}(i)), whence an injective
map $B''\to A^\wedge$. We endow the $\Gamma$-graded ring
$\underline B'':=(B'',\gr_\bullet B'')$ with the topology
induced from $A^\wedge$ via this map, and set
$$
(A,\underline B)^\wedge:=(A^\wedge,\underline B'').
$$

(v)\ \
Notice that, in case $\Gamma\neq 0$, the separation condition
on $\cT$ can be deduced from condition (b) in definition
\ref{def_graded-top} : indeed, if $\gamma,\gamma'$ are two
distinct elements of $\Gamma$, we have
$$
\{0\}^c\subset\gr_\gamma B\cap\gr_{\gamma'}B=0
$$
whence the contention. Thus, the separation condition in
definition \ref{def_graded-top} only serves to rule out
the somewhat trivial case of a non-separated topological
space $A$ endowed with the $\{0\}$-graded $\Z$-algebra
structure.
\end{remark}

\begin{proposition}\label{prop_Cauchy}
Let $(\Gamma,+,0)$ be a monoid, $(A,\underline B)$ a
$\Gamma$-pre-graded structure on the topological ring
$(A,\cT)$. We have :
\begin{enumerate}
\item
The topology $\cT$ is linear.
\item
In the situation of remark {\em\ref{rem_graded-top-algs}(iii)},
the system $(a_S~|~S\subset\Gamma)$, with $S$ ranging over the
filtered set of finite subsets of\/ $\Gamma$, is a Cauchy
net in $B''$, whose unique limit point is $a$.
\item
In the situation of remark {\em\ref{rem_graded-top-algs}(iv)},
the datum $(A,\underline B)^\wedge$ is a $\Gamma$-graded
structure on $(A^\wedge,\cT^\wedge)$.
\item
The inclusion functor $\grTopAlg\to\pregrTopAlg$
admits a left adjoint
$$
\pregrTopAlg\to\grTopAlg
\qquad
(A,\underline B,\Gamma)\mapsto(A,\underline B,\Gamma)^c
$$
that assigns to each topological ring with pre-graded
structure its {\em associated graded structure}.
\end{enumerate}
\end{proposition}
\begin{proof} Fix a fundamental system
$(J_\lambda~|~\lambda\in\Lambda)$ of graded open ideals of $B$.

(i): Since $A$ is separated, the completion map
$(A,\cT)\to(A^\wedge,\cT^\wedge)$ is injective, $B$
is naturally identified with a dense subring of
$A^\wedge$, and the inclusion map $B\to A^\wedge$
extends to an isomorphism between $(A^\wedge,\cT^\wedge)$
and the completion $B^\wedge$ of $B$
(\cite[Ch.II, \S3, n.7, Prop.13]{BouTG}). However,
for every $\lambda\in\Lambda$, let $J^c_\lambda$ be the
topological closure of $J_\lambda$ in $B^\wedge$; then 
$(J^c_\lambda~|~\lambda\in\Lambda)$ is a fundamental
system of open ideals in $B^\wedge$ (remark
\ref{rem_completion-of-topring}(ii)), and lastly,
$(J^c_\lambda\cap A~|~\lambda\in\Lambda)$ is a fundamental
system of open ideals in $A$, whence the contention.

(iii): By remark \ref{rem_general-complete}(iv), the subset
$\gr_\gamma B''$ is closed in the topology of $A^\wedge$,
for every $\gamma\in\Gamma$. Moreover, $B\subset B''$,
so $B''$ is dense in $A^\wedge$. In light of the proof
of (i), it remains only to check that
$J^c_\lambda\cap B''=J''_\lambda$, for every $\lambda\in\Lambda$.
To this aim, say that $x\in B''\cap J^c_\lambda$ for
some such $\lambda$; then there exists a finite subset
$S\subset\Gamma$ such that $x=\sum_{\gamma\in S}x_\gamma$
with $x_\gamma\in\gr_\gamma B''$ for every $\gamma\in S$.
It follows that $\pi^\wedge_\gamma(x)=x_\gamma$ for
$\gamma\in S$ and $\pi^\wedge_\gamma(x)=0$ otherwise.
On the other hand, from claim \ref{cl_image-and-closure} we get
$$
\pi^\wedge_\gamma(x)\in\pi^\wedge_\gamma(J^c_\lambda)
\subset(\pi_\gamma(J_\lambda))^c=\gr_\gamma J^\wedge_\lambda
$$
and the contention follows.

(ii): In the light of the proof of (iii), we have to check
that for every $\lambda\in\Lambda$ there exists a finite
subset $S_J\subset\Gamma$ such that $a_S-a_{S'}\in J''_\lambda$
for every finite subsets $S,S'\subset\Gamma$ containing $S_J$.
However, the proof of (i) shows that the topological closure
$J^c_\lambda$ of $J_\lambda$ in $A$ is an open ideal of $A$;
since $B$ is dense in $A$, it follows that there exists
$b\in B$ such that $a-b\in J^c_\lambda$.
Write $b=\sum_{\gamma\in T}b_\gamma$ for some finite subset
$T\subset\Gamma$; in view of claim \ref{cl_image-and-closure},
we deduce that
$$
\pi^\wedge_\gamma(b-a)=b_\gamma-a_\gamma\in
\pi^\wedge_\gamma(J^c_\lambda)\subset
\pi_\gamma(J)^c=\gr_\gamma J_\lambda^\wedge
\qquad
\text{for every $\gamma\in\Gamma$}
$$
(notation of remark \ref{rem_graded-top-algs}(ii)).
Especially, $a_\gamma\in\gr_\gamma J_\lambda^\wedge$ for every
$\gamma\in\Gamma\setminus T$, and consequently
$\pi^\wedge_\gamma(a_S-a_{S'})\in\gr_\gamma J_\lambda^\wedge$
for every $\gamma\in\Gamma$ and every $S,S'\subset\Gamma$
containing $T$; {\em i.e.} $a_S-a_{S'}\in J''_\lambda$ for
every such $S,S'$, so that $S_J:=T$ will do. By the same
token, $a-a_S=(a-b)+(b-a_S)\in J^c_\lambda$ for every finite
subset $S\subset\Gamma$ containing $T$, so $a$ is the
unique limit of this Cauchy net.

(iv): For every $\gamma\in\Gamma$, let $\gr_\gamma B^c$
be the topological closure of $\gr_\gamma B$ in $A$, and
set $\underline B^c:=\oplus_{\gamma\in\Gamma}\gr_\gamma B^c$;
clearly $\underline B^c\subset\underline B''$, and in
light of (iii), it is clear that the pair $(A,\underline B^c)$
is a topological ring with graded ring structure. Moreover,
if $(f,\phi):(A,\underline B,\Gamma)\to(X,\underline Y,\Gamma')$
is any morphism to a topological ring with graded structure,
then claim \ref{cl_image-and-closure} implies easily that
$f$ maps $\underline B^c$ into $\underline Y$, so we get
the sought left adjoint by setting
$(A,\underline B,\Gamma)^c:=(A,\underline B^c,\Gamma)$.
\end{proof}

\begin{example}\label{ex_was-rem-graded-v}
Let $\Gamma$ be a monoid and $(A,\underline B)$ a topological
ring with $\Gamma$-graded structure. We have the following
elementary operations :

(i)\ \
If $I\subset B$ is any graded ideal, let $I^c$ (resp. $I_A$)
be the topological closure of $I$ in $B$ (resp. in $A$), and
endow $A/I_A$ with the quotient topology arising from the
projection $A\to A/I_A$. Remark \ref{rem_graded-top-algs}(i)
implies that $I^c$ is a graded ideal of $B$ : namely
$\gr_\gamma I^c$ is the topological closure of $\gr_\gamma I$
in $B$, for every $\gamma\in\Gamma$. It follows that
$B/I^c$ is a $\Gamma$-graded subring of $A/I_A$, and
the topology induced by $A/I_A$ on $B/I^c$ agrees with
the quotient topology induced by $B$ (lemma
\ref{lem_no-need-of-dense}(ii)), so the former is linear,
and we get therefore a topological ring with pre-graded
structure $(A/I_A,\underline B/I^c)$. We then may form
the associated graded structure
$$
(A,\underline B)/I:=(A/I_A,\underline B/I^c)^c
$$
(notation of proposition \ref{prop_Cauchy}(iv)). Explicitly,
this is the pair $(A/I_A,\underline C)$, where
$$
\gr_\gamma C:=(\Img(\gr_\gamma B\to A/I_A))^c=
(\gr_\gamma B/I^c)^c
\qquad
\text{for every $\gamma\in\Gamma$}.
$$

(ii)\ \
Keep the situation of (i), and suppose additionally that
$A$ is complete and separated. Then we claim that
$\Img(\gr_\gamma B\to A/I_A)$ is already closed in $A/I_A$,
so the graded subring of $A/I_A$ is just $B/I^c$ in this
case. Indeed, under these assumptions $\gr_\gamma B$ is
also complete and separated, hence the canonical
$\gamma$-projection
$\pi^\wedge_\gamma:A/I_A\to\gr_\gamma C^\wedge$ factors through
the completion map $\gr_\gamma(B/I^c)\to\gr_\gamma C^\wedge$; on
the other hand, $\pi^\wedge_\gamma$ restricts to an injective
map $\gr_\gamma C\to\gr_\gamma C^\wedge$, so there
results an injective map $\gr_\gamma C\to\gr_\gamma(B/I^c)$,
whose restriction to $\gr_\gamma(B/I^c)$ is the identity map.
Thus, $\gr_\gamma C=\gr_\gamma(B/I^c)$, as stated.

(iii)\ \
Lastly, let $\underline C:=(C,\gr_\bullet C)\subset\underline B$
be any closed $\Gamma$-graded subring, and we endow $C$ and the
topological closure $C^c$ of $C$ in $A$ with the topologies
induced from $A$; then $\gr_\gamma C=C\cap\gr_\gamma B$ is closed
in $\gr_\gamma B$, and hence also in $A$, for every $\gamma\in\Gamma$,
so the pair
$$
(A,\underline B)\cap C:=(C^c,\underline C)
$$
is again a topological ring with $\Gamma$-graded structure.
\end{example}

\sset\subsubsection{}\label{subsec_p-powers-Cauchy}
Let $p\geq 2$ be a prime integer, $(A,\underline B)$ any
$\Gamma$-graded topological ring, $a\in A^\wedge$ any element,
and $(a_\gamma~|~\gamma\in\Gamma)$ the sequence attached
to $a$, as in remark \ref{rem_graded-top-algs}(iii).
For every finite subset $S\subset\Gamma$, let $S^*$ be
the set of non-constant mappings $\{1,\dots,p\}\to S$;
choose also a cyclic subgroup $G$ of order $p$ of the
group of permutations of $\{1,\dots,p\}$, and endow
$S^*$ with the right $G$-action given by the rule :
$g(\phi):=\phi\circ g$ for every $g\in G$ and $\phi\in S^*$.
Clearly, for every $\phi\in S^*$, the product
$a_{\phi(i)}\cdots a_{\phi(p)}$ depends only on the class
$\bar\phi$ of $\phi$ in $S^*\!/G$, so we denote it
$a_{\bar\phi}$. Also, notice that every such class $\bar\phi$
has cardinality $p$, since $p$ is a prime; with this notation,
we may then write
$$
(a^p)_S=c_S+p\cdot d_S
\qquad
\text{where $c_S:=\sum_{\gamma\in S}(a_\gamma)^p$
\quad and \quad
$d_S:=\sum_{\bar\phi\in S^*\!/G}a_{\bar\phi}$.}
$$

\begin{lemma}\label{lem__p-power-Cauchy}
With the notation of \eqref{subsec_p-powers-Cauchy}, we have :
\begin{enumerate}
\item
The systems $(c_S~|~S\subset\Gamma)$ and
$(d_S~|~S\subset\Gamma)$ are Cauchy nets in $B''$.
\item
Suppose moreover that the topology of $A$ is coarser than
the $p$-adic topology, and the $p$-Frobenius $\bp_\Gamma$
of\/ $\Gamma$ is injective. Then the following conditions
are equivalent :
\begin{enumerate}
\item
$a$ is topologically nilpotent in $A^\wedge$.
\item
$a_\gamma$ is topologically nilpotent in $B''$ for
every $\gamma\in\Gamma$.
\end{enumerate}
\end{enumerate}
\end{lemma}
\begin{proof}(i): According to proposition \ref{prop_Cauchy}(ii),
for every open graded ideal $J\subset B''$ there exists a
finite subset $S_J\subset\Gamma$ such that
$a_\gamma\in\gr_\gamma J$ for every
$\gamma\in\Gamma\setminus S_J$. Now, if
$S\subset\Gamma$ is any finite subset containing $S_J$,
let $S^*_1$ be the subset of $S^*$ consisting of all mappings
$\{1,\dots,p\}\to S$ whose image is contained in $S_J$,
and set $S_2^*:=S^*\setminus S_1^*$. Clearly
$a_{\bar\phi}\in J$ for every $\bar\phi\in S_2^*\!/G$
(notation of remark \ref{rem_graded-top-algs}(ii)), and
$S^*_1$ is (finite and) independent of $S$, whence
$d_S-d_{S'}\in J$ for every finite subsets
$S,S'\subset\Gamma$ containing $S_J$. The assertion
for $d_S$ follows. Likewise, it is clear that
$c_S-c_{S'}\in J^p$ for every $S,S'$ as in the foregoing,
so the assertion for $c_S$ follows immediately.

(ii): Let $I\subset B$ be any graded open ideal; by assumption
there exists $n\in\N$ such that $p^n\in B$. We consider
the image $\bar a$ of $a$ in the quotient topological
ring with $\Gamma$-graded structure $(A,\underline B)/(I+pB)$,
as in example \ref{ex_was-rem-graded-v}, whose underlying
topological ring is $A/J$, where $J$ is the topological
closure of $I+pB$ in $A$. Let also $I^\wedge$ (resp.
$J^\wedge$) be the topological closure of $I$ (resp.
of $J$) in $A^\wedge$. Then the system of
$\gamma$-projections $(\bar a_\gamma~|~\gamma\in\Gamma)$
of $\bar a$ is the image of $(a_\gamma~|~\gamma\in\Gamma)$.
Now, if $a$ is topologically nilpotent, we may find $k\in\N$
such that $\bar a^{p^k}=0$; since $p\in J$ and since $\bp_\Gamma$
is injective, assertion (i) implies that $(\bar a_\gamma)^{p^k}=0$
for every $\gamma\in\Gamma$, {\em i.e.}
$(a_\gamma)^{p^k}\in J^\wedge$, and therefore
$(a_\gamma)^{p^{kn}}\in(J^\wedge)^n\subset I^\wedge$ (claim
\ref{cl_exchange-pow-clo}). Since $I$ is arbitrary, this
shows that $a_\gamma$ is topologically nilpotent for every
$\gamma\in\Gamma$. Conversely, if the latter condition
holds, for every $\gamma\in\Gamma$ there exists $i\geq n$
such that $\bar a_\gamma^{p^i}=0$. On the other hand, notice
that the topology of $A/J$ is discrete, so $A/J=B/I^c$,
where $I^c$ denotes the topological closure of $I$ in $B$;
especially, there exists a finite set $S\subset\Gamma$ such
that $\bar a_\gamma=0$ for every $\gamma\in\Gamma\setminus S$.
Thus, we may find $j\geq n$ large enough, such that
$\bar a_\gamma^{p^j}=0$ for every $\gamma\in\Gamma$; but
then (i) implies that $\bar a{}^{p^j}=0$, and again, since
$I$ is arbitrary, this implies that $a$ is topologically
nilpotent.
\end{proof}

\begin{example}\label{ex_monoid-alg-graded} 
(i)\ \
Let $\Gamma,\Delta$ be two monoids, and $(A,\underline B)$
a topological ring with $\Delta$-graded structure. Endow
the ring $A[\Gamma]$ with the unique topology such
that the inclusion map $A\to A[\Gamma]$ is adic, and
$B[\Gamma]$ with its natural $\Delta\oplus\Gamma$-graded
structure. Then we claim that the pair
$$
(A,\underline B)[\Gamma]:=(A[\Gamma],B[\Gamma])
$$
is a topological ring with $\Delta\oplus\Gamma$-graded
structure. Indeed, let $(I_\lambda~|~\lambda\in\Lambda)$
be a fundamental system of $\Delta$-graded open ideals
for $B$; then the family
$(I^c_\lambda[\Gamma]~|~\lambda\in\Lambda)$ is a fundamental
system of open ideals for $A[\Gamma]$, and
$I^c_\lambda[\Gamma]\cap B[\Gamma]=I_\lambda[\Gamma]$ for every
$\lambda\in\Lambda$. Moreover, it is easily seen that the
projection $p_\gamma:A[\Gamma]\to A$ :
$\sum_{\gamma'\in\Gamma}\gamma\cdot a_\gamma\mapsto a_\gamma$ is a
continuous map, for every $\gamma\in\Gamma$; now, each direct
summand $\gr_{(\delta,\gamma)}B[\Gamma]$ equals
$p^{-1}_\gamma(\gr_\delta B)\bigcap_{\gamma'\neq\gamma}(\Ker\,p_{\gamma'})$,
so it is a closed subset, whence the contention.

(ii)\ \
Let $\phi:\Gamma\to\Gamma'$ be any morphism of monoids,
$(A,\underline B)$ any topological ring with $\Gamma$-graded
structure. We obtain a topological ring with
$\Gamma'$-pre-graded structure $(A,\underline B{}_{/\Gamma'})$,
and by virtue of proposition \ref{prop_Cauchy}(iv), we may
then form the associated $\Gamma'$-graded structure
$$
(A,\underline B)_{/\Gamma'}:=(A,\underline B{}_{/\Gamma'})^c.
$$

(iii)\ \
In the situation of (ii), notice also that the pair
$(\one_A,\phi)$ yields a morphism of topological rings
with graded structures :
$$
(A,\underline B,\Gamma)\to(A,\underline B,\Gamma)_{/\Gamma'}.
$$
Moreover, suppose that $(A,\underline D,\Gamma')$ is
another $\Gamma'$-graded structure on $A$, such that
$(\one_A,\phi)$ induces a morphism of topological rings
with graded structures
$$
(A,\underline B,\Gamma)_{/\Gamma'}\to
(A,\underline D,\Gamma')
$$
and set $(A,\underline C):=(A,\underline B,\Gamma)_{/\Gamma'}$.
Then we claim that $\underline D=\underline C$.
Indeed, by assumption $\gr_{\gamma'}D$ is closed in
$A$, hence also in $\gr_{\gamma'}C$, for every
$\gamma'\in\Gamma'$; arguing as in remark
\ref{rem_graded-top-algs}(i) we easily see that
$D$ is a closed subset of $C$, whence the contention,
as $D$ is dense in $A$, hence also in $B$.

We deduce that the morphism $(\one_A,\phi)$ enjoys the
following universal property. For every morphism
$(f,\psi):(A,\underline B,\Gamma)\to(A',\underline B',\Gamma')$
such that $\psi$ factors through $\phi$, there exists
a unique morphism of topological rings with $\Gamma'$-graded
structures $(f,\one_{\Gamma'}):(A,\underline B,\Gamma)_{/\Gamma'}
\to(A',\underline B',\Gamma')$ such that
$(f,\psi)=(f,\one_{\Gamma'})\circ(\one_A,\psi)$.
\end{example}

\sset\subsubsection{}
To state the following proposition \ref{prop_general-graded-top},
it is convenient to introduce a special class of monoids; we shall
say that a monoid $(\Gamma,+,0)$ is {\em weakly integral\/} if we
have the following partial cancellation property :
$$
2\cdot\gamma+\delta=2\cdot\gamma+\delta'
\quad\Rightarrow\quad
\gamma+\delta=\gamma+\delta'
\qquad
\text{for every $\gamma,\delta,\delta'\in\Gamma$}.
$$

\begin{lemma} Let $(\Gamma,+,0)$ be any monoid, $n\geq 2$
any integer, and suppose that the $n$-Frobenius endomorphism
of\/ $\Gamma$ is injective. Then $\Gamma$ is weakly integral.
\end{lemma}
\begin{proof} Indeed, let $\gamma,\delta,\delta'\in\Gamma$
be three elements such that
$2\cdot\gamma+\delta=2\cdot\gamma+\delta'$. Then clearly
$n\cdot\gamma+\delta=n\cdot\gamma+\delta'$. A simple
induction on $k$ then shows that
$n\cdot\gamma+k\cdot\delta=n\cdot\gamma+k\cdot\delta'$
for every $k\in\N$. Especially
$n\cdot(\gamma+\delta)=n\cdot(\gamma+\delta')$, whence
$\gamma+\delta=\gamma+\delta'$, as required.
\end{proof}

\begin{lemma}\label{lem_weakly-project}
Let $(\Gamma,+,0)$ be a weakly integral
monoid, $(A,\underline B)$ a topological ring with
$\Gamma$-graded structure, $a\in A$, $\alpha\in\Gamma$
and $f\in\gr_\alpha B$ any three elements. Then we have :
\begin{enumerate}
\item
$\pi^\wedge_{2\alpha+\delta}(f^2a)=f\cdot\pi^\wedge_{\alpha+\delta}(fa)$
for every $\delta\in\Gamma$.
\item
$\pi^\wedge_\gamma(fa)=0$ whenever $\gamma\notin\alpha+\Gamma$.
\end{enumerate}
\end{lemma}
\begin{proof} Since $\pi_\gamma$ is a continuous map for every
$\gamma\in\Gamma$, it suffices to check these identities for
$a\in B$, and then there exists a finite subset
$S\subset\Gamma$ such that $a=a_S$.

(ii): We have $fa=\sum_{\delta\in S}f\cdot\pi^\wedge_\delta(a)$,
and $f\cdot\pi_\delta(a)\in\gr_{\alpha+\delta}B$ for every
$\delta\in S$, whence the assertion.

(i): Clearly $\pi^\wedge_{2\alpha+\delta}(f^2a)=
\sum_{\gamma\in T}a\cdot\pi^\wedge_\gamma(fa)$,
where $T:=\{\gamma\in S~|~\alpha+\gamma=2\alpha+\delta\}$.
But due to (ii), we have $\pi^\wedge_\gamma(fa)=0$ unless
$\gamma\in\alpha+\Gamma$, so we may replace $T$
by the subset of all elements of the form $\alpha+\delta'$
such that $2\alpha+\delta'=2\alpha+\delta$. Since $\Gamma$
is weakly integral, $T=\{\alpha+\delta\}$, whence the
assertion.
\end{proof}

\begin{proposition}\label{prop_general-graded-top}
Let $(\Gamma,+,0)$ be a monoid, $(A,\underline B)$ a
$\Gamma$-graded structure on the topological ring
$(A,\cT)$, and $i_0:\gr_0B\to A$ the inclusion map. We have :
\begin{enumerate}
\item
Suppose that the topology of $A$ is f-adic, complete and
separated. Then :
\begin{enumerate}
\item
There exists a finitely generated graded ideal $J$ of $B$
such that $JA$ is an ideal of adic definition for $A$.
\item
If moreover, $J=(\gr_0J)\cdot B$ or else $\Gamma$ is weakly
integral, $\cT$ induces on $B$ the $J$-adic topology.
\end{enumerate}
\item
Suppose that $i_0$ is c-adic. Then $A$ is c-adic if
and only if the same holds for $\gr_0B$.
\item
Suppose that $A$ is complete and separated, and $i_0$
is c-adic. Then $A$ is f-adic if and only if the same
holds for $\gr_0B$.
\end{enumerate}
\end{proposition}
\begin{proof}(i.a): Let $I$ be any finitely generated ideal
of adic definition of $A$; by assumption, we may find
a graded ideal $I_B$ of $B$ and an integer $n\in\N$ with
$I^n\cap B\subset I_B\subset I$. For every subset $S$ of
$A$ we denote as usual by $S^c$ the topological closure
of $S$ in $A$; since $I^n$ is an open and closed subset
of $A$, we have
$$
I^n=I^n\cap B^c=(I^n\cap B)^c\subset(I^n\cap I_B)^c
$$
so that $I^n=(I^n\cap I_B)+I^{n+1}$. We may then find
a finitely generated graded subideal $J$ of $I_B$
such that $I^n=(I^n\cap J)+I^{n+1}$. On the other
hand, since $A$ is $I$-adically complete and separated,
the ideal $I$ is contained in the Jacobson radical of
$A$ (remark \ref{rem_someth-on-bdd-in-Z-lin}(v)), so
that $I^n=(I^n\cap J)A$, by Nakayama's lemma, and
therefore $I^n\subset JA$; by construction, we have
as well $JA\subset I$, so $J$ is an ideal with the
sought properties.

(i.b): We prove more precisely :

\begin{claim}\label{cl_stimmt}
In the situation of (i), let $J\subset B$ be any graded
ideal, and suppose that either one of the following
conditions holds :
\begin{enumerate}
\alphaenu
\item
$\Gamma$ is integral.
\item
$\Gamma$ is weakly integral and $JA$ is open in $A$.
\item
$J=(\gr_0J)\cdot B$.
\end{enumerate}
Then $J=B\cap JA$ for every $n\in\N$.
\end{claim}
\begin{pfclaim} Say that $x\in B\cap JA$, so  we may
find a mapping $\phi:\{1,\dots,k\}\to\Gamma$ and elements
$a_i\in\gr_{\phi(i)}J$, $y_i\in A$ for $i=1,\dots,k$,
such that $x=\sum_{i=1}^ka_iy_i$, as well as a finite subset
$S\subset\Gamma$ such that $x=x_S$. Suppose first that
$J=(\gr_0J)\cdot B$, in which case we may assume that
$a_i\in\gr_0J$ for $i=1,\dots,k$; we compute
$$
x=\sum_{\gamma\in S}\sum_{i=1}^ka_i\cdot\pi^\wedge_\gamma(y_i)
$$
whence the claim, in this case. Next, suppose that
$\Gamma$ is integral, so that the natural map
$\Gamma\to\Gamma^\gp$ is injective; we may then replace
$\Gamma$ by $\Gamma^\gp$, and assume from start that
$\Gamma$ is an abelian group. We compute
$$
x=\sum_{\gamma\in S}\sum_{i=1}^ka_i\cdot\pi^\wedge_{\gamma-\phi(i)}(y_i)
$$
so the claim follows also in this case. Lastly, suppose
that condition (b) holds; then we may find homogeneous
elements $f_1,\dots,f_r\in J$ such that $\sum_{i=1}^rf_iA$
is an open ideal of $A$, and the same holds therefore for
the ideal $I:=\sum_{i=1}^rf^2_iA$. Next, we may find $y'_i\in B$
such that $y_i-y'_i\in I$ for $i=1,\dots,k$, and clearly
it suffices to check that $x':=x-\sum_{i=1}^ky'_ia_i\in J$.
We may then replace $x$ by $x'$, and assume as well that
$x\in I$; in this case, lemma \ref{lem_weakly-project}
says that $\pi_\gamma(x)\in\sum_{i=1}^rf_iB$ for every
$\gamma\in S$, whence the contention.
\end{pfclaim}

(ii): Due to lemma \ref{lem_still-c-adic}(i) and proposition
\ref{prop_Cauchy}(iii), we may replace $(A,\underline B)$
by $(A,\underline B)^\wedge$, after which we may assume
that $A$ is complete and separated.
We know already that if the topology of $\gr_0B$
is c-adic, the same holds for $\cT$ (lemma
\ref{lem_f-adics}(i.b)). Thus, suppose that $\cT$ is
c-adic, and let $I$ be any ideal of c-adic definition
of $A$; by assumption, there exists an open ideal
$I_0\subset\gr_0B$ with $(I_0A)^c\subset I$. Now, let
$J\subset\gr_0B$ be any other open ideal; we know that
there exists $m\in\N$ such that $(I^m)^c\subset(JA)^c$,
and taking into account claim \ref{cl_exchange-pow-clo}
we deduce that $(I_0^mA)^c\subset(JA)^c$. Combining
with claim \ref{cl_image-and-closure} we get
$$
\pi^\wedge_0(I_0^mA)^c\subset\pi^\wedge_0(JA)^c.
$$
But notice that, since $A$ is complete,
$\pi^\wedge_0(I_0^mA)=I_0^m$ and $\pi^\wedge_0(JA)=J$, so
finally $(I^m_0)^c\subset J$, whence the assertion.

(iii): Suppose first that $\gr_0B$ is f-adic (in which
case, it is also adic, by proposition \ref{prop_f-adics}(iii)).
Then lemmata \ref{lem_f-adics}(i.b) and \ref{lem_5.3.8}(ii.a)
imply that the same holds for $A$. Conversely, suppose
that $A$ is f-adic (and adic), and let $I\subset A$ be
any finitely generated ideal of adic definition. By
assumption, there exists an open ideal $J\subset\gr_0B$
such that $I^n\subset (JA)^c\subset I$ for some $n\in\N$.
Since $I^n$ is both open and closed in $A$, we deduce
$$
I^n=I^n\cap(JA)^c=(I^n\cap JA)^c=I^{n+1}+(JA\cap I^n).
$$
It follows that there exists a finitely generated graded
subideal $K\subset J$ such that $I^n\subset KA+I^{n+1}$, and
so $I^n\subset KA\subset I$, by virtue of Nakayama's
lemma and remark \ref{rem_someth-on-bdd-in-Z-lin}(v).
Then $(K^rA~|~r\in\N)$ is a fundamental system of open
ideals in $A$, and therefore $(K^rA\cap\gr_0B~|~r\in\N)$
is a fundamental system of open ideals in $\gr_0B$. But
from claims \ref{cl_exchange-pow-clo} and
\ref{cl_image-and-closure} we get
$$
K^r\subset K^rA\cap\gr_0B\subset\pi^\wedge_0(K^rA)\subset
\pi^\wedge_0(K^rB)^c=K^r
$$
so $K^rA\cap\gr_0B=K^r$ for every $r\in\N$,
and finally $\gr_0B$ is f-adic.
\end{proof}

\subsection{Homological algebra for topological modules}
\label{sec_Topolog-Hom-Alg}
Let $A$ be any topological ring. The category of topological
$A$-modules and continuous $A$-linear maps, denoted
$$
A\tdu\TopMod
$$
is additive and with representable kernels and cokernels but,
generally, not abelian; indeed, if $f:M\to N$ is a continuous
map of topological $A$-modules, $\Coker(\Ker\,f\to M)$ is not
necessarily isomorphic to $\Ker(N\to\Coker\,f)$, since the
quotient topology (induced from $M$) on $\Img\,f$ may be finer
than the subspace topology (induced from $N$). It is therefore
useful to introduce the following :

\begin{definition}\label{def_strict-map}
Let $A$ be a topological ring, $f:M\to N$ a morphism of
topological $A$-modules. We say that $f$ is {\em strict},
if the natural map
$$
\Coker(\Ker\,f\to M)\to\Ker(N\to\Coker\,f)
$$
is an isomorphism of topological $A$-modules (where these
kernels and cokernels are formed in the additive category
$A\tdu\TopMod$ : see remark \ref{rem_additive-cat}(v)).
\end{definition}

The category $A\tdu\TopMod$ is {\em exact\/} in the sense of
\cite{Qu}; namely, the admissible monomorphisms
(resp. epimorphisms) are the continuous injections (resp.
surjections) $f:M\to N$ that are strict, in the sense of
definition \ref{def_strict-map}, {\em i.e.} that induce
homeomorphisms $M\isom f(M)$ (resp. $M/\Ker\,f\isom N$),
where $f(M)$ (resp. $M/\Ker\,f$) is endowed with the subspace
(resp. quotient) topology induced from $N$ (resp. from $M$).
An admissible epimorphism is also called a {\em quotient map}
of topological $A$-modules. Correspondingly there is a well
defined class of {\em admissible short exact sequences\/} of
topological $A$-modules. The following lemma exhibits a useful
class of admissible short exact sequences.

\begin{lemma}\label{lem_crit-admissible}
Consider an inverse system, with surjective transition maps
$$
(\underline E_n~|~n\in\N)\quad :\quad
0\to(M'_n~|~n\in\N)\to(M_n~|~n\in\N)\to(M''_n~|~n\in\N)\to 0
$$
of short exact sequences of discrete $A$-modules. Then
the induced complex of inverse limits :
$$
\lim_{n\in\N}\underline E_n\quad :\quad
0\to M':=\lim_{n\in\N}M'_n\to M:=\lim_{n\in\N}M_n\to
M'':=\lim_{n\in\N}M''_n\to 0
$$
is an admissible short exact sequence of topological $A$-modules.
\end{lemma}
\begin{proof} For every pair of integers $i,j\in\N$ with $i\leq j$,
let $\phi_{ji}:M_j\to M_i$ be the transition map in the inverse
system $(M_n~|~n\in\N)$, and define likewise $\phi'_{ji}$ and
$\phi''_{ji}$; the assumption means that all these maps are
onto. Set $K_{ji}:=\Ker\,\phi_{ji}$ and define likewise $K'_{ji}$,
$K''_{ji}$. By the snake lemma we deduce, for every $i\in\N$, an
inverse system of short exact sequences :
$$
0\to(K'_{ji}~|~j\geq i)\to(K_{ji}~|~j\geq i)\to
(K''_{ji}~|~j\geq j)\to 0
$$
where again, all the transition maps are surjective. However,
by definition, the decreasing family of submodules
$(K_i:=\lim_{j\geq i}K_{ji}~|~i\in\N)$ is a fundamental system
of open neighborhoods of $0\in M$ (and likewise
one defines the topology on the other two inverse limits).
It follows already that the surjection $M\to M''$ is a
quotient map of topological $A$-modules. To conclude it
suffices to remark :

\begin{claim} The natural map $K'_i:=\lim_{j\geq i}K'_{ji}\to K_i$
induces an identification :
$$
K'_i=K_i\cap M'\qquad\text{for every $i\in\N$}.
$$
\end{claim}
\begin{pfclaim}[] Indeed, for every $j\geq i$ we have a left
exact sequence :
$$
\underline F_j\quad :\quad (0\to K'_{ji}\to M'_j\oplus K_{ji}
\stackrel{\alpha}{\longrightarrow} M_j)
$$
where $\alpha(m',m)=m'-m$ for every $m'\in M'_j$ and $m\in K_{ji}$.
Therefore $\lim_{j\geq i}\underline F_j$ is the analogous left exact
sequence $0\to K'_i\to M'\oplus K_i\to M$, whence the claim.
\end{pfclaim}
\end{proof}

\begin{remark}\label{rem_admiss-complexes}
Let $(C^\bullet,d^\bullet)$ be a complex of topological $A$-modules
with continuous differentials. Then, for every $i\in\Z$, the
cohomology $H^iC^\bullet$ inherits a well defined topology : namely,
we can take either
$$
\Coker(d^{i-1}:C^{i-1}\to\Ker\,d^i)
\qquad\text{or}\qquad
\Ker(d^i:\Coker\,d^{i-1}\to C^{i+1})
$$
(where these kernels and cokernels are formed in the
additive category $A\tdu\TopMod$) and by lemma
\ref{lem_no-need-of-dense}(i), these two topologies
coincide (regardless of whether the differentials are
strict or otherwise). Notice also that for every $i\in\Z$
the following conditions are equivalent :
\begin{enumerate}
\alphaenu
\item
The differential $d^i$ factors through an open map
$C^i\to\Ker\,d^{i+1}$.
\item
The differential $d^i$ is strict and the topology
of $H^{i+1}C^\bullet$ is discrete.
\end{enumerate}
\end{remark}

\sset\subsubsection{}\label{subsec_adm-sho-ex-seq-cplx}
Consider an {\em admissible short exact sequence}
of complexes of topological $A$-modules
$$
0\to(C^\bullet_1,d^\bullet_1)\xrightarrow{\ f^\bullet\ }
(C^\bullet_2,d^\bullet_2)\xrightarrow{\ g^\bullet\ }
(C^\bullet_3,d^\bullet_3)\to 0
$$
{\em i.e.} a double complex of $A\tdu\TopMod$ whose rows
$0\to C^k_1\to C^k_2\to C^k_3\to 0$ are admissible short
exact sequences for every $k\in\Z$. For every $i\in\Z$
and for $j=1,2,3$, denote by $\delta^i_j:C^i\to\Ker\,d^{i+1}_j$
the continuous map induced by the differential
$d^i_j:C^i_j\to C^{i+1}_j$ of $C^\bullet_j$.

\begin{proposition}\label{prop_when-deltas-open}
In the situation of \eqref{subsec_adm-sho-ex-seq-cplx}, fix
also $i\in\Z$. We have :
\begin{enumerate}
\item
If $\delta^i_2$ and $\delta^{i-1}_3$ are open maps, the same
holds for $\delta^i_1$.
\item
If $\delta^i_1$ and $\delta^i_3$ are open maps, the same holds
for $\delta^i_2$.
\item
If $\delta^i_1$ and $\delta^{i-1}_2$ are open maps, the same
holds for $\delta^{i-1}_3$.
\end{enumerate}
\end{proposition}
\begin{proof}(i): Endow $D:=g^{i+1}(\Ker\,d_2^{i+1})$ with the
quotient topology induced by the restriction
$h:\Ker\,d_2^{i+1}\to D$ of $g^{i+1}$, and notice that
$\Img\,d^i_3\subset D$; we claim that $d^i_3$ induces
a continuous and open map $C^i_3\to D$. Indeed, since $g^i$
is an admissible epimorphism, it suffices to check that
the composition $d^i_3\circ g^i:C^i_2\to D$ is continuous
and open; however, the latter equals $h\circ\delta^i_2$,
whence the assertion. Moreover, $f^{i+1}$ restricts to an
admissible monomorphism $\Ker\,d^{i+1}_1\to\Ker\,d^{i+1}_2$.
Summing up, we may replace $C^{i+1}_j$ with $\Ker\,d^{i+1}_j$
for $j=1,2$, and $C^{i+1}_3$ with $D$, and assume from start
that $C^{i+2}_j=0$ for $j=1,2,3$. Then $d^i_j:C^i_j\to C^{i+1}_j$
is open for $j=2,3$, and it remains to check that $d^i_1$ is
open as well.

Next, set $\bar C{}^i_1:=\Coker\,d^{i-1}_1$ and
$\bar C{}^i_2:=\Coker\,d^{i-1}_2\circ f^{i-1}$ (endowed with
the topologies induced by $C^i_1$ and respectively $C^i_2$).
Clearly $d^i_j$ factors through a continuous map
$\bar d{}^i_j:\bar C{}^i_j\to C^{i+1}_j$ for $j=1,2$, and $f^i$
(resp. $g^i$) induces an injective (resp. surjective)
continuous map $\bar f{}^i:\bar C{}^i_1\to\bar C{}^i_2$
(resp. $\bar g{}^i:\bar C{}^i_2\to C^i_3$). Moreover, by
lemma \ref{lem_no-need-of-dense}(i) the topology of
$\bar C{}^i_1$ is induced by the topology of $\bar C{}^i_2$
via $\bar f{}^i$, {\em i.e.} $\bar f{}^i$ is an admissible
monomorphism. Likewise, $\bar g{}^i$ is an admissible
epimorphism. We may then replace $C^i_j$ by $\bar C{}^i_j$
for $j=1,2$, and assume from start $C^{i-1}_1=0$, in which
case we obtain a commutative diagram of continuous maps
\set\begin{equation}\label{eq_cont-admiss-maps}
{\diagram 0 \ar[r] & 0 \ar[r] \ar[d] &
C^{i-1}_2 \ar[r]^-{g^{i-1}} \ar[d]_{d_2^{i-1}} &
C^{i-1}_3 \ar[d]^{d_3^{i-1}} \ar[r] & 0 \\
0 \ar[r] & C^i_1 \ar[r]^-{f^i} \ar[d]_{d^i_1} &
C^i_2 \ar[r]^-{g^i} \ar[d]_{d^i_2} &
C^i_3 \ar[r] \ar[d]^{d^i_3} & 0 \\
0 \ar[r] & C^{i+1}_1 \ar[r]^-{f^{i+1}} &
C^{i+1}_2 \ar[r]^-{g^{i+1}} & C^{i+1}_3 \ar[r] & 0
\enddiagram}
\end{equation}
whose rows are admissible short exact sequences, and where
$d^i_2$ and $\delta_3^{i-1}$ are still open.
Consider the diagram of topological $A$-modules
\set\begin{equation}\label{eq_new-corner}
{\diagram
\Ker\,(d^i_3\circ g^i) \ar[r]^-{\bar g{}^i} \ar[d]_{\bar d{}^i_2} &
\Ker\,d^i_3 \ar[d] \\
\Ker\,g^{i+1} \ar[r] & 0
\enddiagram}
\end{equation}
deduced from the right bottom square of \eqref{eq_cont-admiss-maps}.
In view of claim \ref{cl_trivial-open} we deduce that $\bar g{}^i$
and $\bar d{}^i_2$ are open maps. Moreover, $f^i$ and $f^{i+1}$
factor through admissible monomorphisms
$C^i_1\to\Ker\,(d^i_3\circ g^i)$ and $C^{i+1}_1\to\Ker\,g^{i+1}$.
Likewise, $d^{i-1}_2$ and $d^{i-1}_3$ factor through continuous maps
$C^{i-1}_2\to\Ker\,(d^i_3\circ g^i)$ and $C^{i-1}_3\to\Ker\,d^i_3$.
Thus, we may replace the bottom right square of
\eqref{eq_cont-admiss-maps} with \eqref{eq_new-corner}, and
assume from start that $C^{i+1}_3=0$ as well. In this situation,
$g^{i-1}$ and $f^{i+1}$ are both isomorphisms of topological
$A$-modules, and $d^{i-1}_3$ is an open map, so we are reduced
to showing :

\begin{claim}\label{cl_cross-diagram}
Consider a diagram of topological $A$-modules
$$
\xymatrix{ & & C_4 \ar[d]^h \\
0 \ar[r] & C_1 \ar[r]^-f & C_2 \ar[r]^-g \ar[d]^k
& C_3 \ar[r] & 0 \\
& & C_5
}$$
such that :
\begin{enumerate}
\alphaenu
\item
The horizontal row is an admissible short exact sequence.
\item
$g\circ h$ and $k$ are open maps and $k\circ h=0$. 
\end{enumerate}
Then $k\circ f$ is an open map.
\end{claim}
\begin{pfclaim} After pulling back the short exact admissible
sequence of $(f,g)$ along the map $g\circ h$, we obtain
the short exact sequence of topological modules
$$
0\to C_1\xrightarrow{\ f'\ } C_2\times_{C_3}C_4
\xrightarrow{\ g'\ } C_4\to 0.
$$
Notice that $g'$ is still an admissible epimorphism, by
claim \ref{cl_trivial-open}. Likewise, $f'$ is an admissible
monomorphism, since its composition with the projection
$\pi:C_2\times_{C_3}C_4\to C_2$ equals $f$, which is an
admissible monomorphism by assumption. Moreover, since
$g\circ h$ is open, the same holds for $\pi$, again by
claim \ref{cl_trivial-open}, and therefore also for
$k':=k\circ\pi:C_2\times_{C_3}C_4\to C_5$. Furthermore,
the pair $(h,\one_{C_4})$ defines a continuous section
$s:C_4\to C_2\times_{C_3}C_4$ for $g'$, whence a continuous
isomorphism of $A$-modules
$$
\omega:C_1\times C_4\isom C_2\times_{C_3}C_4
\qquad
(x,y)\mapsto f'(x)+s(y)=(f(x)+h(y),y)
$$
and we notice that $\omega$ is an isomorphism of
topological $A$-modules : indeed, let us endow $f(C_1)$
with the topology induced by the inclusion into $C_2$,
so that $f$ factors through an isomorphism $u:C_1\isom f(C_1)$
of topological $A$-modules; then the inverse of $\omega$
is the continuous map
$$
C_2\times_{C_3}C_4\to C_1\times C_4
\qquad
(a,b)\mapsto(u^{-1}(a-h(b)),b)
$$
whence the assertion. Summing up, we deduce that the continuous
map $\psi:=k\circ\pi\circ\omega:C_1\times C_4\to C_5$ is open,
and notice that $\psi(x,y)=k\circ f(x)$ for every
$(x,y)\in C_1\times C_4$, since $k\circ h=0$. Lastly, let
$U\subset C_1$ be any open subset; then $U\times C_4$ is
open in $C_1\times C_4$, and $\psi(U\times C_4)=k\circ f(U)$
is open in $C_5$, whence the claim.
\end{pfclaim}

(ii): Define the $A$-module $D$ and the $A$-linear map
$h:\Ker\,d^{i+1}_2\to D$ as in (i), and endow $D$ with the
topology induced by the inclusion into $C^{i+1}_3$. Then $h$
is continuous, and it is also an open map, since the same
holds for $h\circ\delta^i_2=\delta^i_3\circ g^i:C^i_2\to D$
(lemma \ref{lem_descent-of-openness}).
We may then argue as in the proof of (i), to reduce to
the case where $C^{i+2}_j=0$ for $j=1,2,3$, in which case
$d^i_1$ and $d^i_3$ are both open maps, and we must show
that the same holds for $d^i_2$.

Next, we pull back the admissible short exact sequence
$(f^i,g^i)$ (resp. $(f^{i+1},g^{i+1})$) along the map $g^i$
(resp. along the map $d^i_3\circ g^i$) to get the
commutative diagram of topological $A$-modules:
\set\begin{equation}\label{eq_many-pullbacks}
{\diagram 0 \ar[r] & C^i_1 \ar[r]^-{f'^i} \ar[d]_{d^i_1} &
C^i_2\times_{C^i_3}C^i_2 \ar[r]^-{g'^i} \ar[d]_{d'^i_2} &
C^i_2 \ddouble \ar[r] & 0 \\
0 \ar[r] & C^{i+1}_1 \ar[r]^-{f'^{i+1}} \ddouble &
C^{i+1}_2\times_{C^{i+1}_3}C^i_2 \ar[r]^-{g'^{i+1}} \ar[d]_{\pi^{i+1}} &
C^i_2 \ar[d]^{d^i_3\circ g^i} \ar[r] & 0 \\
0 \ar[r] & C^{i+1}_1 \ar[r]^-{f^{i+1}} &
C^{i+1}_2 \ar[r]^-{g^{i+1}} & C^{i+1}_3 \ar[r] & 0.
\enddiagram}
\end{equation}
Arguing as in the proof of claim \ref{cl_cross-diagram}
we easily see that the rows of \eqref{eq_many-pullbacks}
are still admissible short exact sequences. Moreover, the
pair $(\one_{C^i_2},\one_{C^i_2})$ (resp. $(d^i_2,\one_{C^i_2})$)
gives a section $s^i$ (resp. $s^{i+1}$) for $g'^i$ (resp.
for $g'^{i+1}$). Furthermore, $\pi^{i+1}$ is an open map,
by claim \ref{cl_trivial-open}, and if 
$\pi^i:C^i_2\times_{C^{i+1}_3}C^i_2\to C^i_2$ denotes the
projection, we have
$$
d^i_2\circ\pi^i=\pi^{i+1}\circ d'^i_2.
$$
Thus, in order to show that $d^i_2$ is open, it suffices
to check that the same holds for $d'^i_2$ (lemma
\ref{lem_descent-of-openness}). Summing up, we may replace
the exact rows $(f^i,g^i)$ and $(f^{i+1},g^{i+1})$ by
$(f'^i,g'^i)$ and $(f'^{i+1},g'^{i+1})$, and assume also from
start that $C^i_3=C^{i+1}_3$ and $d^i_3=\one_{C^i_3}$, and
moreover, that $g^i$ and $g^{i+1}$ admit continuous sections
$s^i$ and respectively $s^{i+1}$, such that $d^i_2\circ s^i=s^{i+1}$.
Arguing as in the proof of claim \ref{cl_cross-diagram}
we deduce isomorphisms of topological $A$-modules
$$
\omega^k:C^k_1\times C^i_3\isom C^k_2
\qquad
(x,y)\mapsto f^k(x)+s^k(y)
\qquad
\text{for $k=i,i+1$}
$$
fitting into a commutative diagram
$$
\xymatrix{
C^i_1\times C^i_3 \ar[rr]^-{\omega^i} \ar[d]_{d^i_1\times\one_{C^i_3}}
& & C^i_2 \ar[d]^{d^i_2} \\
C^{i+1}_1\times C^i_3 \ar[rr]^-{\omega^{i+1}} & & C^{i+1}_2
}$$
and since $d^i_1$ is an open map, the same holds for
$d^i_1\times\one_{C^i_3}$, whence the assertion.

(iii): Arguing as in the proof of (i), we reduce to the case
where $C^{i+2}_j=0$ for every $j=1,2,3$, in which case $d^i_1$
is an open map. Then, we define $\bar C{}^i_j$ for $j=1,2$ as
in the proof of (i), so that $d^i_j$ factors through a
continuous map $\bar d{}^i_j:\bar C{}^i_j\to C^{i+1}_j$ for
$j=1,2$. Notice that $\bar d{}^i_1$ is an open map, and
$\delta^{i-1}_2$ induces an open map
$\Coker\,f^{i-1}\to\Ker\,\bar d{}^i_2$. We may then replace
$C^i_j$ by $\bar C{}^i_j$ for $j=1,2$, and assume from
start that $C^{i-1}_1=0$, in which case we obtain a
commutative diagram \eqref{eq_cont-admiss-maps} whose rows
are admissible short exact sequences, and where $d^i_1$ and
$\delta^{i-1}_2$ are open maps. Next, we may argue again as
in the proof of (i), to reduce to the case where $C^{i+1}_3=0$
as well, in which case we have to show that $d^{i-1}_3$ is
an open map. To this aim, we consider the cartesian diagram
$$
\xymatrix{
C^i_1\times_{C^{i+1}_2}C^i_2 \ar[rr]^-{\pi^i_2} \ar[d]_{\pi^i_1} & &
C^i_2 \ar[d]^{d^i_2} \\
C^i_1 \ar[rr]^-{f^{i+1}\circ d^1_1} & & C^{i+1}_2
}$$
and we notice that the pair $(\one_{C^i_1},f^i)$ defines a
continuous section $s:C^i_1\to C^i_1\times_{C^{i+1}_2}C^i_2$ for
$\pi^i_1$. Likewise, the pair consisting of the zero map
$\Ker\,d^i_2\to C^i_1$ and the inclusion map $\Ker\,d^i_2\to C^i_2$
defines a continuous map
$t:\Ker\,d^i_2\to C^i_1\times_{C^{i+1}_2}C^i_2$, whence an isomorphism
of topological $A$-modules
$$
\omega:C^i_1\times\Ker\,d^i_2\isom C^i_1\times_{C^{i+1}_2}C^i_2
\qquad
(x,y)\mapsto s(x)+t(y)=(x,f^i(x)+y).
$$
Furthermore, under the current assumptions, $f^{i+1}$ is an
isomorphism of topological $A$-modules, so $f^{i+1}\circ d^i_1$
is an open map, and then the same holds for $\pi^i_2$, according
to claim \ref{cl_trivial-open}. Lastly, notice the commutative
diagram
$$
\xymatrix{ C^i_1\times C^{i-1}_2 \ar[r]^-{\pi^{i-1}_2}
\ar[d]_{\one_{C^i_1}\times\delta^{i-1}_2}
& C^{i-1}_2 \ar[rr]^-{g^{i-1}} & & C^{i-1}_3 \ar[d]^{d^{i-1}_3} \\
C^i_1\times\Ker\,d^i_2 \ar[r]^-\omega &
C^i_1\times_{C^{i+1}_3}C^i_2 \ar[r]^-{\pi_2^i} &
C^i_2 \ar[r]^-{g^i} & C^i_3
}$$
where $\pi^{i-1}_2$ is the projection. Since $\delta^{i-1}_2$ is
an open map, the same holds for $\one_{C^i_1}\times\delta^{i-1}_2$,
and then also for $d^{i-1}_3\circ g^{i-1}\circ\pi^{i-1}_2$. By
lemma \ref{lem_descent-of-openness}, we conclude that $d^{i-1}_3$
is open as well, as required.
\end{proof}

\sset\subsubsection{}\label{eq_strict-cont-double}
Consider an object $C^{\bullet\bullet}$ of $\sC_2(A\tdu\TopMod)$,
{\em i.e.} a double complex of topological $A$-modules whose
horizontal and vertical differentials
$$
d_h^{p,q}:C^{p,q}\to C^{p+1,q}
\qquad
d_v^{p,q}:C^{p,q}\to C^{p,q+1}
$$
are continuous maps, and therefore induce continuous maps
$$
\delta^{p,q}_h:C^{p,q}\to\Ker\,d_h^{p+1,q}
\qquad
\delta^{p,q}_v:C^{p,q}\to\Ker\,d_v^{p,q+1}
\qquad
\text{for every $p,q\in\Z$}.
$$
With this notation, we have :

\begin{corollary}\label{cor_double-strict-continuous}
In the situation of \eqref{eq_strict-cont-double},
suppose moreover that :
\begin{enumerate}
\alphaenu
\item
$C^{p,q}=0$ whenever $p<0$ or $q<0$.
\item
The maps $\delta^{p,q}_v$ are open for every $p,q\in\Z$.
\item
The maps $\delta^{p,q}_h$ are open for every $p\in\Z$ and
every $q>0$.
\end{enumerate}
Then the maps $\delta^{p,0}_h$ are open for every $p\in\Z$.
\end{corollary}
\begin{proof} To ease notation, set
$$
B^{pq}_v:=\Img\,d^{p,q-1}_v
\qquad
Z^{pq}_v:=\Ker\,d^{p,q}_v
\qquad
H^{p,q}_v:=\Coker\,\delta^{p,q}_h
\qquad
\text{for every $p,q\in\Z$}
$$
and denote by $d_B^{p,q}:B^{p,q}_v\to B^{p+1,q}$  and
$d^{p,q}_Z:Z^{p,q}_v\to Z^{p+1,q}$ the restrictions of
$d^{p,q}_h$, for every $p,q\in\Z$. Then
$(B^{\bullet,q}_v,d_B^{\bullet,q})$ and $(Z^{\bullet,q}_v,d_Z^{\bullet,q})$
are complexes of topological $A$-modules, and we get
admissible short exact sequences of complexes
$$
\begin{aligned}
& 0\to(Z^{\bullet,q}_v,d_Z^{\bullet,q})\to(C^{\bullet,q},d_h^{\bullet,q})\to
(B^{\bullet,q+1}_v,d_B^{\bullet,q+1})\to 0 \\
& 0\to(B^{\bullet,q}_v,d_B^{\bullet,q})\to(Z^{\bullet,q}_v,d_Z^{\bullet,q})\to
H^{\bullet,q}_v\to 0.
\end{aligned}
$$
Let also
$$
\delta_B^{p,q}:B^{p,q}_v\to\Ker d_B^{p+1,q}
\qquad
\delta_Z^{p,q}:Z^{p,q}_v\to\Ker d_Z^{p+1,q}
$$
be the continuous maps induced by $d_B^{p,q}$ and respectively
$d_Z^{p,q}$, for every $p,q\in\Z$.

\begin{claim}\label{cl_bootstrap}
For every $p,q$ the following holds :
\begin{itemize}
\item[(d)]
If $\delta_Z^{p,q}$ is an open map, the same holds for
$\delta_B^{p,q}$.
\item[(e)]
If $\delta_B^{p-1,q+1}$ is an open map, the same holds for
$\delta_Z^{p,q}$, provided $q>0$.
\end{itemize}
\end{claim}
\begin{pfclaim} Indeed, $H^{\bullet,q}_v$ is complex of discrete
topological $A$-modules, for every $q\in\Z$, due to assumption
(b), so its differentials $d_H^{\bullet,q}$ trivially induce open
maps $H^{p,q}_v\to\Ker\,d_H^{p+1,q}$ for every $q$; then
(d) follows from proposition \ref{prop_when-deltas-open}(i).
Likewise, (e) follows from assumption (b) and proposition
\ref{prop_when-deltas-open}(i).
\end{pfclaim}

We now show that, for every $n\in\Z$, the maps $\delta_B^{p,q}$
are open for every $p,q\in\Z$ such that $p+q=n$ and $q>0$.
We argue by induction on $p$. For $p<-1$, assumption (a) says
that $B_v^{p,q}=B_v^{p+1,q}=0$, so the assertion trivially holds.
Suppose then that $p\geq-1$, and that we already know that
$\delta_B^{p-1,q+1}$ is open; from claim \ref{cl_bootstrap}
we deduce that $\delta_B^{p,q}$ is also open, as stated.
Especially, we see that $\delta_B^{p,1}$ is an open map, for
every $p\in\Z$.
Lastly, notice that (a) and (b) also imply that the topology
of $Z_v^{p,0}$ is discrete for every $p\in\Z$; then we apply
proposition \ref{prop_when-deltas-open}(ii) to the admissible
short exact sequence
$0\to Z_v^{\bullet,0}\to C^{\bullet,0}\to B_v^{\bullet,1}$ to conclude.
\end{proof}
 
\sset\subsubsection{}\label{subsec_contract-continuously}
Let $A$ be any topological ring, $(C^\bullet,d^\bullet)$ any
complex of topological $A$-modules, $\phi^\bullet:C^\bullet\to C^\bullet$
an endomorphism of $C^\bullet$ and $h^\bullet$ a chain homotopy
from $\one_{C^\bullet}$ to $\phi^\bullet$ in the category
$\sC(A\tdu\TopMod)$, so that $\phi^i:C^i\to C^i$ and
$h^i:C^i\to C^{i-1}$ are continuous $A$-linear maps, for every
$i\in\Z$.

\begin{lemma}\label{lem_contract-cont}
In the situation of \eqref{subsec_contract-continuously},
fix $i\in\Z$ and suppose that $\phi^i=0$. Then we have :
\begin{enumerate}
\item
$d^{i-1}$ induces an open and surjective map
$\delta^i:C^{i-1}\to\Ker\,d^i$.
\item
$\Ker\,d^i$ is a direct summand of the topological $A$-module
$C^i$ ({\em i.e.} the inclusion map $\Ker\,d^i\to C^i$ admits
a continuous left inverse).
\end{enumerate}
\end{lemma}
\begin{proof} The assumption means that
$$
d^{i-1}\circ h^i+h^{i+1}\circ d^i=\one_{C^i}.
$$
Now, set $e^i:=d^{i-1}\circ h^i$; we compute
$$
(\one_{C^i}-e^i)\circ e^i=h^{i+1}\circ d^i\circ d^{i-1}\circ h^i=0
$$
and therefore $e^i$ is an idempotent continuous endomorphism
of $C^i$. According to example \ref{ex_idemp-in-additive-cat},
it follows that $e^i$ and $\one_{C^i}-e^i$ induce a natural
isomorphism of topological $A$-modules
\set\begin{equation}\label{eq_cont-decompose}
\Ker\,(e^i)\oplus\Ker\,(\one_{C^i}-e^i)\isom C^i.
\end{equation}
Notice as well that $d^i\circ e^i=0$, whence
$d^i=d^i\circ(\one_{C^i}-e^i)$, and we deduce easily that
$$
\Ker\,(\one_{C^i}-e^i)=\Ker\,d^i
$$
whence (ii). By the same token, we see that $e^i$ restricts
to the identity map on $\Ker\,d^i$, {\em i.e.} $\delta^i$
admits a continuous right inverse, given by the restriction
$\Ker\,d^i\to C^{i-1}$ of $h^i$. Taking into account
lemma \ref{lem_descent-of-openness}, we then get also (i).
\end{proof}

\sset\subsubsection{Extensions of topological rings}
\label{subsec_exaltops}
Let $A$ be any topological ring whose topology is linear;
we shall consider for any $A$-algebra $C$, the category
$$
\sExal_A(C)
$$
whose objects are all short exact sequences of $A$-modules
\set\begin{equation}\label{eq_exalt}
\Sigma
\quad :\quad
0\to M\to E\xrightarrow{\ \psi\ } C\to 0
\end{equation}
such that $E$ is an $A$-algebra (with $A$-module structure
given by the structure map $A\to E$), and $\psi$ is a map
of $A$-algebras. The morphisms in $\sExal_A(C)$ are the
commutative ladders of $A$-modules whose central vertical
arrow $g$ is a map of $A$-algebras :
\set\begin{equation}\label{eq_morph-in-exalt}
{\diagram
0 \ar[r] & M' \ar[r] \ar[d] & E' \ar[r] \ar[d]^g &
C \ar[r] \ddouble & 0 \\
0 \ar[r] & M'' \ar[r] & E'' \ar[r] & C \ar[r] & 0.
\enddiagram}
\end{equation}
For any topological $A$-algebra $C$ whose topology is linear,
we shall consider the category
$$
\sExaltop_A(C)
$$
whose objects are the short exact sequences of $A$-modules
\eqref{eq_exalt} such that $E$ is a topological $A$-algebra
whose topology is linear (and again, with $A$-module structure
given by the structure map $A\to E$), and $\psi$ is a
continuous and open map of topological $A$-algebras, whose
kernel $M$ is a discrete topological space, for the topology
induced from $E$. The morphisms in $\sExaltop_A(C)$ are the
commutative ladders \eqref{eq_morph-in-exalt} such that $g$
is a continuous map of topological $A$-algebras.

\sset\subsubsection{}\label{subsec_ext-top.algebras}
Now, let $A$ be as in \eqref{subsec_exaltops}, $C$ any
topological $A$-algebra whose topology is linear, $\phi:A\to C$
the structure morphism, and $(I_\lambda~|~\lambda\in\Lambda)$
(resp. $(J_{\lambda'}~|~\lambda'\in\Lambda')$) a collection of
open ideals of $A$ (resp. of $C$), which is a fundamental
system of open neighborhoods of $0$. Let
$\Lambda''\subset\Lambda\times\Lambda'$ be the subset of all
$(\lambda,\lambda')$ such that $\phi(I_\lambda)\subset J_{\lambda'}$;
the set $\Lambda''$ is partially ordered, by declaring that
$(\lambda,\lambda')\leq(\mu,\mu')$ if and only if
$I_\mu\subset I_\lambda$ and $J_{\mu'}\subset J_{\lambda'}$.
Set $A_\lambda:=A/I_\lambda$ for every $\lambda\in\Lambda$
(resp. $C_{\lambda'}:=C/J_{\lambda'}$ for every $\lambda'\in\Lambda'$);
for $(\lambda,\lambda'),(\mu,\mu')\in\Lambda''$ with
$(\lambda,\lambda')\leq(\mu,\mu')$, the surjection
$\pi_{\mu'\lambda'}:C_{\mu'}\to C_{\lambda'}$ induces a functor
$$
\sExal_{A_\lambda}(C_{\lambda'})\to\sExal_{A_\mu}(C_{\mu'})
\qquad
\Sigma\mapsto\Sigma*\pi_{\mu'\lambda'}
$$
(see \cite[\S 2.5.5]{Ga-Ra}) and clearly the rule
$(\lambda,\lambda')\mapsto\sExal_{A_\lambda}(C_{\lambda'})$
yields a pseudo-functor
$$
\sE:(\Lambda'',\leq)\to\bCat.
$$
Moreover, we have a pseudo-cocone :
\set\begin{equation}\label{eq_pseudo-cocco}
\sE\Rightarrow\sExaltop_A(C)
\end{equation}
defined as follows. To any $(\lambda,\lambda')\in\Lambda''$ and any
object $\Sigma_{\lambda,\lambda'}\!:\!
0\!\!\to\!\!M\!\!\to\!E_{\lambda'}\!\!\to\!C_{\lambda'}\!\!\to\!0$
of $\sExal_{A_\lambda}(C_{\lambda'})$, one assigns the extension
\eqref{eq_exalt} obtained by pulling back $\Sigma_{\lambda,\lambda'}$
along the projection $\pi_{\lambda'}:C\to C_{\lambda'}$. Let
$\beta:E\to E_{\lambda'}$ be the induced projection; we endow
$E$ with the linear topology defined by the fundamental system
of all open ideals of the form $\psi^{-1}J\cap\beta^{-1}J'$ where
$J$ (resp. $J'$) ranges over the set of open ideals of $C$ (resp.
over the set of all ideals of $E_{\lambda'}$). With this topology,
it is easily seen that both $\psi$ and the induced structure map
$A\to E$ are continuous ring homomorphisms. Moreover, if
$J\subset J_{\lambda'}$, then
$\psi^{-1}J\cap\beta^{-1}J'=J\times(J'\cap M)$, from which it
follows that $\psi$ is an open map. Furthermore,
$M\cap\beta^{-1}0=0$, which shows that $M$ is discrete, for the
topology induced by the inclusion map $M\to E$. Summing up, we
have attached to $\Sigma_{\lambda,\lambda'}$ a well defined
object $\Sigma_{\lambda,\lambda'}*\pi_{\lambda'}$ of
$\sExaltop_A(C)$, and it is easily seen that the rule
$\Sigma_{\lambda,\lambda'}\mapsto\Sigma_{\lambda,\lambda'}*\pi_{\lambda'}$
is functorial in $\Sigma_{\lambda,\lambda'}$, and for
$(\lambda,\lambda)\leq(\mu,\mu')$ in $\Lambda''$, there is
a natural isomorphism in $\sExaltop_A(C)$ :
$$
\Sigma_{\lambda,\lambda'}*\pi_{\lambda'}\isom
(\Sigma_{\lambda,\lambda'}*\pi_{\mu'\lambda'})*\pi_{\mu'}
$$
(details left to the reader).

\begin{lemma}\label{lem_exaltops}
The pseudo-cocone \eqref{eq_pseudo-cocco} induces an equivalence
of categories :
$$
\beta:\Pscolim{\Lambda''}\sE\isom\sExaltop_A(C).
$$
\end{lemma}
\begin{proof} Let $\Sigma$ as in \eqref{eq_exalt} be any object
of $\sExaltop_A(C)$. By assumption, there exists an open ideal
$J\subset E$ such that $M\cap J=0$. Since $\psi$ is an open
map, there exists $\lambda'\in\Lambda'$ such that
$J_{\lambda'}\subset\psi(J)$, and after replacing $J$ by
$J\cap\psi^{-1}J_{\lambda'}$, we may assume that
$\psi(J)=J_{\lambda'}$. Likewise, if $\phi_E:A\to E$ is the
structure morphism, there exists $\lambda\in\Lambda$ such that
$I_\lambda\subset\phi_E^{-1}J$, and it follows that the induced
extension
$$
\Sigma_{\lambda,\lambda'}:0\to M\to E/J\to C_{\lambda'}\to 0
$$
is an object of $\sExal_{A_\lambda}(C_{\lambda'})$. We notice :

\begin{claim}\label{cl_natural-top.iso}
There exists a natural isomorphism
$\Sigma\isom\Sigma_{\lambda,\lambda'}*\pi_{\lambda'}$ in
$\sExaltop_A(C)$.
\end{claim}
\begin{pfclaim}(i): By construction, $\psi$ and the projection
$\pi_J:E\to E/J$ define a unique morphism
$\gamma:E\to(E/J)\times_{C_{\lambda'}}C$ of $A$-algebras,
restricting to the identity map on $M$ (which is an ideal
in both of these $A$-algebras). It is clear that $\gamma$
is an isomorphism, and therefore it yields a natural isomorphism
$\Sigma\isom\Sigma_{\lambda,\lambda'}*\pi_{\lambda'}$ in
$\sExal_A(C)$. It remains to check that $\gamma$ is continuous
and open. For the continuity, it suffices to remark that, for
every ideal $I\subset E/J$ and every $\mu'\in\Lambda'$,
the ideals $\gamma^{-1}(I\times_{C_{\lambda'}}C)=\pi_J^{-1}I$
and $\gamma^{-1}(E/J\times_{C_{\lambda'}}J_{\mu'})=\psi^{-1}J_{\mu'}$
are open in $E$, which is obvious, since $J$ is an open ideal
and $\psi$ is continuous. Lastly, let $I\subset E$ be any open
ideal such that $I\subset J\cap\psi^{-1}J_{\lambda'}$; since
$\psi$ is an open map, it is easily seen that
$\gamma(I)=0\times\psi(I)$ is an open ideal of
$(E/J)\times_{C_{\lambda'}}C$, so $\gamma$ is open.
\end{pfclaim}

From claim \ref{cl_natural-top.iso} we see already that $\beta$
is essentially surjective. It also follows easily that $\beta$
is full. Indeed, consider any morphism $s:\Sigma'\to\Sigma''$ of
$\sExaltop_A(C)$ as in \eqref{eq_morph-in-exalt}, and pick an
open ideal $J''\subset E''$ with $J\cap M''=0$; set $J':=g^{-1}J''$,
and notice that $J'\cap M'=0$.
Moreover, if the image of $J''$ in $C$ equals $J_{\lambda'}$ for
some $\lambda'\in\Lambda'$, then clearly the same holds for the
image of $J'$ in $C$. Therefore, in this case the foregoing
construction yields objects $\Sigma'_{\lambda,\lambda'}$ and
$\Sigma''_{\lambda,\lambda'}$ of $\sExal_{A_\lambda}(C_{\lambda'})$
(for a suitable $\lambda\in\Lambda$), whose middle terms are
respectively $E'/J'$ and $E''/J''$, and $s$ descends to a morphism
$s_{\lambda'}:\Sigma'_{\lambda,\lambda'}\to\Sigma''_{\lambda,\lambda'}$,
whose middle term is the map $g_{\lambda'}:E'/J'\to E''/J''$ induced
by $g$. By inspecting the proof of claim \ref{cl_natural-top.iso},
we deduce a commutative diagram
$$
\xymatrix{ E' \ar[r]^-\sim \ar[d]_g &
(E'/J')\times_{C_{\lambda'}}C \ar[d]^{g_{\lambda'}\times_{C_{\lambda'}}C} \\
E'' \ar[r]^-\sim & (E''/J'')\times_{C_{\lambda'}}C
}$$
whose horizontal arrows are the maps that define the isomorphisms
$\Sigma'\isom\Sigma'_{\lambda,\lambda'}*\pi_{\lambda'}$ and
$\Sigma''\isom\Sigma''_{\lambda,\lambda'}*\pi_{\lambda'}$ in
$\sExaltop_A(C)$. It follows easily that
$s_{\lambda}*\pi_{\lambda'}=s$, whence the assertion.
Lastly, the faithfulness of $\beta$ is immediate, since the
projections $\pi_{\lambda'}$ are surjective maps.
\end{proof}

\sset\subsubsection{}
Let $A$ be as in \ref{subsec_exaltops}, and $C$ any
$A$-algebra (resp. any topological $A$-algebra); we denote by
$$
\snilExal_A(C)
\qquad
\text{(\ resp.\ \ $\snilExaltop_A(C)$\ )}
$$
the full subcategory of $\sExal_A(C)$ (resp. of $\sExaltop_A(C)$)
whose objects are the {\em nilpotent extensions\/} of $C$,
{\em i.e.} those extensions \eqref{eq_exalt}, where $M$
is a nilpotent ideal of $E$. Moreover, in the situation of
\eqref{subsec_ext-top.algebras}, clearly $\sE$ restricts to a
pseudo-functor
$$
\mathsf{nilE}:(\Lambda'',\leq)\to\bCat
\qquad
(\lambda,\lambda')\mapsto\snilExal_{A_\lambda}(C_{\lambda'}).
$$
Also, \eqref{eq_pseudo-cocco} restricts to a pseudo-cocone on
$\mathsf{nilE}$, and lemma \ref{lem_exaltops} immediately implies
an equivalence of categories
\set\begin{equation}\label{eq_pseudo-nil-cocco}
\Pscolim{\Lambda''}\mathsf{nilE}\isom\snilExaltop_A(C).
\end{equation}

\begin{proposition}\label{prop_exalt}
Let $A$ and $C$ be as in \eqref{subsec_ext-top.algebras}.
The following holds :
\begin{enumerate}
\item
Suppose that $J^2_{\lambda'}$ is open in $C$, for every
$\lambda'\in\Lambda'$. Then the forgetful functor
\set\begin{equation}\label{eq_forget-top-ext}
\snilExaltop_A(C)\to\snilExal_A(C)
\end{equation}
is fully faithful.
\item
Suppose additionally, that :
\begin{enumerate}
\item
$C$ is noetherian, and $I\subset C$ is an ideal such that
the topology of\/ $C$ is $I$-adic.
\item
$I^2_\lambda$ is open in $A$, for every $\lambda\in\Lambda$.
\end{enumerate}
Then the essential image of \eqref{eq_forget-top-ext} is the
(full) subcategory of all nilpotent extensions \eqref{eq_exalt}
such that the $C$-module $M/M^2$ is annihilated by a power of $I$.
\end{enumerate}
\end{proposition}
\begin{proof}(i): The functor is obviously faithful, and in light
of \eqref{eq_pseudo-nil-cocco}, we come down to the following
situation. We have a commutative ladder of extensions of $A$-algebras
$$
\xymatrix{
0 \ar[r] & M \ar[r] \ar[d] & E \ar[r]^-\psi \ar[d]^g &
C \ar[r] \ar[d]^{\pi_{\lambda'}} & 0 \\
0 \ar[r] & N \ar[r] & E' \ar[r] & C_{\lambda'} \ar[r] & 0
}$$
for some $\lambda'\in\Lambda'$, whose top (resp. bottom) row
is an object of $\snilExaltop_A(C)$ (resp. of
$\snilExal_{A_\lambda}(C_{\lambda'})$, for some $\lambda\in\Lambda$),
and we need to show that the kernel of $g$ contains an open ideal.
To this aim, we remark :

\begin{claim}\label{cl_square-open}
For any open ideal $J\subset E$, the ideal $J^2$ is open as well.
\end{claim}
\begin{pfclaim} We may assume that $J\cap M=0$. We have
$I:=\psi(J^2)=\psi(J)^2$, so the assumption in (i) say that
$I$ is open in $C$, and therefore $\psi^{-1}I=J^2\oplus M$
is open in $E$, so finally $J\cap\psi^{-1}I=J^2$ is open in
$E$, as stated.
\end{pfclaim}

Now, set $J:=\psi^{-1}J_{\lambda'}$; then $J$ is an open
ideal of $E$, and clearly $g(I)\subset N$. Say that $N^k=0$;
then $g(I^k)=0$, and $I^k$ is an open ideal of $E$, by claim
\ref{cl_square-open}, whence the contention.

(ii): It is easily seen that, for every extension \eqref{eq_exalt}
in the essential image of \eqref{eq_forget-top-ext}, we must
have $I^k(M/M^2)=0$ for every sufficiently large $k\geq 0$ (details
left to the reader). Conversely, consider a nilpotent extension
$\Sigma$ as in \eqref{eq_exalt}, with $I^k(M/M^2)=0$ and $M^t=0$
for some $k,t\in\N$. Pick a finite system $\bff:=(f_1,\dots,f_r)$
of elements of $E$ whose images in $B$ form a system of generators
for $I^k$; notice that the ideal $(\bff^n)B$ is open in $B$, and
$(\bff^n)M=0$ for every $n\geq t$. There follows an inverse system
of exact sequences (notation of remark \ref{rem_koszul-alg}(ii))
$$
H_1(\bff^n,B)\to M\xrightarrow{\ \beta_n\ }E/(\bff^n)E\to
B/(\bff^n)B\to 0
\qquad
\text{for every $n\geq t$}
$$
with transition maps induced by the morphisms $\bphi_\bff$
of \eqref{subsec_def-bphis}. As $B$ is noetherian, lemma
\ref{lem_Hartsho} and remark \ref{rem_Hartsho} imply that
the system $(H_1(\bff^n,B)~|~n\in\N)$ is essentially zero.
Therefore, the same holds for the inverse system
$(\Ker\,\beta_n~|~n>0)$. However, the transition map
$\Ker\,\beta_{n+1}\to\Ker\,\beta_n$ is obviously injective
for every $n\geq t$, so we conclude that $\beta_n$ is injective
for some sufficiently large integer $n$. For such $n$, we
obtain a nilpotent extension
\set\begin{equation}\label{eq_final-exaltop}
0\to M\to E/(\bff^n)E\to B/(\bff^n)B\to 0.
\end{equation}
Lastly, let $\phi:A\to B$ be the structure map, and pick
$\lambda\in\Lambda$ such that $I_\lambda\subset\phi^{-1}I^k$;
it is easily seen that $I^t_\lambda M=0$ and
$\phi(I^s_\lambda)\subset(\bff^n)B$ for $s\in\N$ large
enough. Therefore, some sufficiently large power of $I_\lambda$
annihilates $E/(\bff^n)E$; under our assumption (b), such
power of $I_\lambda$ contains another open ideal $I_\mu$,
so \eqref{eq_final-exaltop} is an object of
$\snilExal_{A_\mu}(B/(\bff^n)B)$ whose image in
$\snilExaltop_A(B)$ agrees with $\Sigma$.
\end{proof}

\begin{proposition}\label{prop_fsmooth-cot}
Let $A$ be a topological ring (whose topology is linear), $B$
a noetherian $A$-algebra, $I\subset B$ an ideal, and suppose that :
\begin{enumerate}
\alphaenu
\item
The structure map $A\to B$ is continuous for the $I$-adic
topology on $B$.
\item
For every open ideal $J\subset A$, the ideal $J^2$ is also open.
\end{enumerate}
Then the following conditions are equivalent :
\begin{enumerate}
\alphaenu
\addenu\addenu
\item
$B$ (with its $I$-adic topology) is a formally smooth $A$-algebra.
\item
$\Omega^1_{B/A}\otimes_BB/I$ is a projective $B/I$-module, and
$H_1(\L_{B/A}\otimes_BB/I)=0$.
\end{enumerate}
\end{proposition}
\begin{proof} More generally, let $A$ and $C$ be as in
\eqref{subsec_ext-top.algebras}, and $M$ a {\em discrete\/}
$C$-module, {\em i.e.} a $C$-module annihilated by an open
ideal; we denote by $\Exaltop_A(C,M)$ the $C$-module of
square zero topological $A$-algebra extensions of $C$ by $M$.
Likewise, for every $(\lambda,\lambda')\in\Lambda''$, and
every $C_{\lambda'}$-module 
$M$, let $\Exal_{A_\lambda}(C_{\lambda'},M)$ be the
$C_{\lambda'}$-module of isomorphism classes of square zero
$A_\lambda$-algebra extensions of $C_{\lambda'}$ by $M$.
For $(\lambda,\lambda')\leq(\mu,\mu')$, we get a natural map of
$C_{\mu'}$-modules
\set\begin{equation}\label{eq_exal-grave}
\Exal_{A_\lambda}(C_{\lambda'},M)\to\Exal_{A_\mu}(C_{\mu'},M)
\qquad
\text{for every $C_{\lambda'}$-module $M$}
\end{equation}
and \eqref{eq_pseudo-nil-cocco} implies a natural isomorphism :
$$
\colim_{(\lambda,\lambda')\in\Lambda''}
\Exal_{A_\lambda}(C_{\lambda'},M)\isom\Exaltop_A(C,M)
$$
which is well defined for every discrete $C$-module $M$.
By inspecting the definitions, and taking into account
\cite[Ch.0, Prop.19.4.3]{EGAIV}, it is easily seen that $C$ is a
formally smooth $A$-algebra if and only if $\Exaltop_A(C,M)=0$
for every such discrete $C$-module $M$. We also have a natural
$C$-linear map :
\set\begin{equation}\label{eq_forget-topology}
\Exaltop_A(C,M)\to\Exal_A(C,M)
\qquad
\text{for every discrete $C$-module $M$}
\end{equation}
and proposition \ref{prop_exalt}(i) shows that
\eqref{eq_forget-topology} is an injective map, provided
$J_{\lambda'}^2$ is an open ideal of $C$, for every
$\lambda'\in\Lambda'$. Moreover, for $C:=B$ (with its
$I$-adic topology), the map \eqref{eq_forget-topology} is
an isomorphism, by virtue of proposition \ref{prop_exalt}(ii).

Now, recall the natural isomorphism of $C_\lambda$-modules
(\cite[Ch.III, Th.1.2.3]{Il})
$$
\Exal_A(B,M)\isom\Ext^1_B(\L_{B/A},M)
\qquad
\text{for every $B$-module $M$}.
$$
Combining with the foregoing, we deduce a natural $B$-linear
isomorphism
$$
\Exaltop_A(B,M)\isom\Ext^1_B(\L_{B/A},M)
$$
for every $B$-module $M$ annihilated by some ideal $I^k$. Summing
up, $B$ is a formally smooth $A$-algebra if and only if
$\Ext^1_B(\L_{B/A},M)$ vanishes for every discrete $B$-module
$M$. Set $B_0:=B/I$; by considering the $I$-adic filtration on
a given $M$, a standard argument shows that the latter condition
holds if and only if it holds for every $B_0$-module $M$, and in
turns, this is equivalent to the vanishing of
$\Ext^1_{B_0}(\L_{B/A}\otimes_BB_0,M)$ for every $B_0$-module
$M$. This last condition is equivalent to (d), as stated.
\end{proof}

We conclude this section with a list of topological
corollaries of the conditions
$\mathrm{(a)}_\bff-\mathrm{(f)}_\bff$ that were introduced
in \eqref{subsec_badabum}.

\sset\subsubsection{}\label{subsec_Urquart}
Let $A$ be a ring, $\bff:=(f_1,\dots,f_r)$ a finite
sequence of elements of $A$ that generate an ideal
$I\subset A$. Let also $(C_\bullet,d_\bullet)$ be any
bounded above complex of flat $A$-modules such that,
for every $j\in\Z$ the annihilator ideal of
$H_jC_\bullet$ contains a power of $I$. Denote also by
$A^\wedge$ the $I$-adic completion of $A$, and by
$C^\wedge_\bullet$ the complex whose term (resp.
differential) in degree $j$ is the $I$-adic completion
of $C_j$ (resp. of $d_j$), for every $j\in\Z$.
Lastly, set $Q_\bullet:=\Cone(A[0]\to A^\wedge[0])$,
the cone of the morphism given by the completion map,
and let $\bff^\wedge$ be the image in $A^\wedge$ of the
sequence $\bff$.

\begin{corollary}\label{cor_cond-a-and-compl}
In the situation of \eqref{subsec_Urquart}, the following
holds :
\begin{enumerate}
\item
If $A$ satisfies condition $\mathrm{(a)_\bff}$ of
\eqref{subsec_badabum}, the natural maps
$$
C_\bullet\to A^\wedge\otimes_AC_\bullet\to C^\wedge_\bullet
$$
are quasi-isomorphisms.
\item
The following conditions are equivalent :
\begin{enumerate}
\item
$A$ satisfies condition $\mathrm{(a)_\bff}$.
\item
$A^\wedge$ satisfies condition $\mathrm{(a)_{\bff^\wedge}}$
and $Q_\bullet\derotimes_AA/I[0]$ is acyclic.
\end{enumerate}
\end{enumerate}
\end{corollary}
\begin{proof}(i): According to \cite[Th.3.5.8]{We}, for
every $i\in\Z$ we have a natural short exact sequence
$$
0\to L:=\lim_{n\in\N}{}^{\!1}H_{i+1}(C_\bullet\otimes_AA/I^n)\to
H_i(C^\wedge)\xrightarrow{\ \omega\ }
\lim_{n\in\N}H_i(C_\bullet\otimes_AA/I^n)\to 0.
$$
On the other hand, the complex $C_\bullet$ fulfills
condition $\mathrm{(f)_\bff}$ of \eqref{subsec_badabum},
so the kernel and cokernel of the natural morphism of
inverse systems
$(H_i(C_\bullet)~|~n\in\N)\to(H_i(C_\bullet\otimes_AA/I^n)~|~n\in\N)$
are both essentially zero, hence the induced map
$$
H_i(C_\bullet)\to\lim_{n\in\N}H_i(C_\bullet\otimes_AA/I^n)
$$
is an isomorphism for every $i\in\Z$, and $L=0$
(\cite[Prop.3.5.7]{We}). Hence, $\omega$ is an isomorphism,
and it is easily seen that the resulting isomorphism
$H_i(C^\wedge)\isom H_i(C)$ is the inverse of the map
induced by the completion map $C_\bullet\to C^\wedge_\bullet$,
hence the latter is a quasi-isomorphism. Next we remark :

\begin{claim}\label{cl_two-equiv-Tor-vanish}
Let $B$ be any ring, $\bg$ a finite sequence of elements
of $B$ that generates an ideal $J$, and $K_\bullet$ a bounded
above complex of $B$-modules. The following conditions are
equivalent :
\begin{enumerate}
\alphaenu
\item
$K_\bullet\derotimes_BB/J[0]$ is acyclic.
\item
$\bK_\bullet(\bg,K_\bullet)$ is acyclic.
\end{enumerate}
\end{claim}

\begin{pfclaim} (a)$\Rightarrow$(b): For every $q\in\N$
we consider the spectral sequences
$$
\begin{aligned}
E^2_{pq}:=&\, H_p(K_\bullet\derotimes_BH_q\bK_\bullet(\bg)[0])
\Rightarrow H_{p+q}(\bg,K_\bullet) \\
F(q)^2_{ij}:=&\,
\Tor_i^{B/J}(H_j(K_\bullet\derotimes_BB/J[0]),H_q\bK_\bullet(\bg))
\Rightarrow E^2_{i+j,q} & \text{for every $q\in\Z$}.
\end{aligned}
$$
Assumption (a) implies that $F(q)^2_{ij}=0$ for every $i,j,q\in\Z$,
so that $E^2_{pq}=0$ for every $p\in\Z$, whence (b). Conversely,
if (b) holds, we show -- by induction on $n$ -- that
$H_n(K_\bullet\derotimes_BB/J[0])=0$ for every $n\in\Z$.
First, since $K_\bullet$ is bounded above, we may find $n_0\in\Z$
such that the assertion holds for every $n\leq n_0$. Next,
let $k>n_0$, and suppose that the assertion is already known
for every $n<k$. It follows that $F(q)^2_{ij}=0$ whenever
$j<k$, and clearly the same holds also whenever $i<0$;
thus $F(q)^2_{0k}=E^2_{kq}$ for every $q\in\Z$. By the same
token, we deduce that $E^2_{pq}=0$ for every $p<k$, and
clearly the same holds whenever $q<0$. Summing up, it
follows that $E^2_{k0}=H_k(\bg,K_\bullet)=0$; however,
$F(0)^2_{0k}=H_k(K_\bullet\derotimes_BB/J[0])$, whence the
contention.
\end{pfclaim}
\begin{claim}\label{cl_more-will-vanish}
Suppose that $A$ satisfies condition $\mathrm{(a)_\bff}$
of \eqref{subsec_badabum}, let $k\in\N$ be any integer,
and $N$ any $A/I^k$-module. The completion map
$A\to A^\wedge$ induces an isomorphism
$$
N\isom A^\wedge\otimes_AN
\qquad\text{and}\qquad
\Tor^A_i(A^\wedge,N)=0
\qquad
\text{for every $i>0$}.
$$
\end{claim}
\begin{pfclaim} Notice that, since $\bK_j(\bff)$ is a
free $A$-module of finite rank for every $j\in\Z$, the
natural map
$A^\wedge\otimes_A\bK_\bullet(\bff)\to\bK_\bullet(\bff)^\wedge$
is an isomorphism. But we have already shown that the
completion map $\bK_\bullet(\bff)\to\bK_\bullet(\bff)^\wedge$
is a quasi-isomorphism, hence $\bK_\bullet(\bff,Q_\bullet)$ is
acyclic, and in light of claim \ref{cl_two-equiv-Tor-vanish}
we see that the same holds for $Q_\bullet\derotimes_AA/I[0]$.
Let $M$ be any $A/I$-module; from the change of ring spectral
sequence
$$
\Tor^{A/I}_i(H_j(Q_\bullet\otimes_AA/I[0]),M)\Rightarrow
H_{i+j}(Q_\bullet\derotimes_AM[0])
$$
we deduce that $Q_\bullet\derotimes_AM[0]$ is also
acyclic. Lastly, by considering, for $t=0,\dots,k-1$
the distinguished triangles
$$
Q_\bullet\derotimes_AI^tN/I^{t+1}N[0]\to
Q_\bullet\derotimes_AN/I^{t+1}N[0]\to
Q_\bullet\derotimes_AN/I^tN[0]\to
Q_\bullet\derotimes_AI^tN/I^{t+1}N[1]
$$
a simple induction shows that $Q_\bullet\derotimes_AN[0]$
is acyclic as well; the latter assertion is equivalent
to the claim.
\end{pfclaim}

Now, let us consider the spectral sequence
$$
E^2_{pq}:=\Tor^A_p(A^\wedge,H_qC_\bullet)\Rightarrow
H_{p+q}(A^\wedge[0]\derotimes_AC_\bullet)
$$
and notice that $A^\wedge[0]\derotimes_AC_\bullet=
A^\wedge\otimes_AC_\bullet$, since $C_\bullet$ is a complex
of flat $A$-modules. By assumption, for every $q\in\Z$
there exists $k\in\N$ such that $I^k\cdot H_qC_\bullet=0$;
taking into account claim \ref{cl_more-will-vanish} we
conclude that $E^2_{pq}=0$ whenever $p>0$, and
$E^2_{0q}=H_qC_\bullet$ for every $q\in\Z$; it is then
easily seen that the resulting isomorphism
$H_qC_\bullet\isom H_q(A^\wedge\otimes_AC_\bullet)$ is the
map induced by the completion map $A\to A^\wedge$.

(ii.a)$\Rightarrow$(ii.b): We have already remarked
in the proof of (i) that if $A$ satisfies condition
$\mathrm{(a)_\bff}$, then $Q_\bullet\derotimes_AA/I[0]$
is acyclic. Next, define $A_r$ and $\beta_\bff:A_r\to A$
as in \eqref{subsec_badabum}, and consider the spectral
sequence
$$
E(n)^2_{ij}:=H_i(\Tor^{A_r}_j(A_r/I^n_r,A)\derotimes_AQ_\bullet)
\Rightarrow
H_{i+j}(A_r/I^n_r\derotimes_{A_r}Q_\bullet)
\qquad
\text{for every $n\in\N$}.
$$
In light of claim \ref{cl_more-will-vanish} we see that
$E(n)^2_{ij}=0$ for every $i\in\Z$ and $j,n\in\N$, so
$A_r/I^n_r\derotimes_{A_r}Q_\bullet$ is acyclic, and therefore
the completion map $A\to A^\wedge$ induces quasi-isomorphisms
\set\begin{equation}\label{eq_used-twice}
A_r/I^n_r\otimes_{A_r}A\isom A_r/I^n_r\otimes_{A_r}A^\wedge
\qquad
\text{for every $n\in\N$}
\end{equation}
so $A^\wedge$ satisfies $\mathrm{(a)_\bff}$.

(ii.b)$\Rightarrow$(ii.a): Arguing as in the proof of
claim \ref{cl_more-will-vanish} we see that
$N\derotimes_AQ_\bullet$ is acyclic for every $k\in\N$
and every $A/I^k$-module $N$. Then we find again
$E(n)^2_{ij}=0$ for every $i\in\Z$ and every $j,n\in\N$,
as well as the quasi-isomorphisms \eqref{eq_used-twice},
and the assertion follows.
\end{proof}

\begin{corollary}\label{cor_go-to-quotient}
In the situation of \eqref{subsec_go-to-quotient}, set
$Z:=\Spec\,A/I$, and suppose that $U:=\Spec\,A\setminus Z$
is an affine scheme. Then $\bar A$ satisfies condition
$\mathrm{(c)^{un}_{\bar\bff}}$ of \eqref{subsec_badabum}.
\end{corollary}
\begin{proof} Endow $A$ with its $I$-adic topology, and
$A_U:=\cO_U(U)$ with the f-adic topology $\cT_U$ provided
by proposition \ref{prop_top-on-opens-fadic-case}(i).
Clearly, the restriction map $\rho_U:A\to A_U$ factors
through $\bar A$, so the latter is an open subring of
$A_U$. Now, by definition there is no prime ideal of $U$
that contains all the elements $\bar f_1,\dots,\bar f_r$;
since $U$ is affine by assumption, it follows that there
exist $g_1,\dots,g_r\in A_U$ such that
$\sum_{i=1}^rg_i\cdot\bar f_i=1$ in $A_U$. Next, since
$\bar A$ is a ring of definition of $A_U$, and $\bar I:=I/J$
an ideal of adic definition, it follows that there exists
$k\in\N$ such that $g_i\cdot\bar I{}^k\subset\bar A$ for
$i=1,\dots,r$. It follows easily that the scalar
multiplication by $g_i$ induces an $\bar A$-linear map
$\bar I{}^{n+k}\to\bar I{}^n$ for every $n\in\N$, and the
latter in turn induces a map of complexes
$$
\psi_{i,n}:\bK_\bullet(\bar\bff,\bar I^{n+k})\to
\bK_\bullet(\bar\bff,\bar I{}^n)
\qquad
\text{for every $n\in\N$ and $i=1,\dots,r$}.
$$
On the other hand, scalar multiplication by $\bar f_i$
induces a map of complexes
$$
\phi_{i,n}:\bK_\bullet(\bar\bff,\bar I{}^n)\to
\bK_\bullet(\bar\bff,\bar I{}^n)
\qquad
\text{for every $n\in\N$ and $i=1,\dots,r$}
$$
that equals the zero morphism in the homotopy category
(see lemma \ref{lem_koszul-vanish}(i)). Summing up,
we conclude that $\sum_{i=1}^r\phi_{i,n}\circ\psi_{i,n}$
is the zero morphism $\bK_\bullet(\bar\bff,\bar I{}^{n+k})
\to\bK_\bullet(\bar\bff,\bar I{}^n)$ in the homotopy
category, and on the other hand it coincides with the
inclusion map. The assertion follows.
\end{proof}

\begin{remark}\label{rem_Ann-and-completion}
(i)\ \
Let $A$ be a ring, $J\subset A$ an ideal of finite type; endow
$A$ with its $J$-adic topology, denote by $A^\wedge$ the completion
of $A$, and set $X_A:=\Spec\,A$, $X_{A^\wedge}:=\Spec\,A^\wedge$.
Let also $U\subset X_A$ be an affine open subset containing the
analytic locus of $X_A$ (see definition \ref{def_deja-vu}), and
set $U^\wedge:=U\times_{X_A}X_{A^\wedge}$. Moreover, let
$$
A_U:=\cO_{X_A}(U)
\quad
A_{U^\wedge}:=\cO_{X_{A^\wedge}}(U^\wedge)
\quad
\tau_A:=\Ker(A\to A_U)
\quad
\tau_{A^\wedge}:=\Ker(A^\wedge\to A_{U^\wedge})
$$
and denote by $(\tau_AA^\wedge)^c$ the topological closure of
$\tau_AA^\wedge$ in $A^\wedge$. Then we have
$$
\tau_{A^\wedge}\subset(\tau_AA^\wedge)^c.
$$
Indeed, endow $A_U$ with the f-adic topology characterized by
proposition \ref{prop_top-on-opens-fadic-case}(i), and let
$A^\wedge_U$ be the completion of $A_U$; then $A/\tau_A$ is
a ring of definition of $A_U$, hence its completion
$(A/\tau_A)^\wedge=A^\wedge/(\tau_AA^\wedge)^c$ is a ring of
definition of $A^\wedge_U$, and the natural map
$A^\wedge/(\tau_AA^\wedge)^c\otimes_AA_U\to A^\wedge_U$ is an
isomorphism (propositions \ref{prop_replaces-Mat-Th.8.1}(iii)
and \ref{prop_complete-f-adic}(i,ii)). On the other hand,
$A_{U^\wedge}=A^\wedge\otimes_AA_U$, so we get a commutative
diagram of rings :
$$
\xymatrix{
A^\wedge \ar[r] \ar[d] & A^\wedge/(\tau_AA^\wedge)^c \ar[d] \\
A_{U^\wedge} \ar[r] & A^\wedge_U 
}$$
whose right vertical arrow is injective, whence the assertion.

(ii)\ \
In the situation of (i), say furthermore that
$X_A\setminus U=\Spec\,A/I$, for some ideal $I\subset A$
generated by a finite sequence $\bff:=(f_1,\dots,f_r)$ of
elements of $A$, and suppose that $\tau_A\cdot I^k=0$ for some
$k\in\N$. Then $\tau_A\cap I^n=0$ for some $n\in\N$. Indeed,
by virtue of corollary \ref{cor_go-to-quotient} and proposition
\ref{prop_a-un-and-quots}, the ring $A$ satisfies condition
$\mathrm{(c)^{un}_{\bff}}$ of \eqref{subsec_badabum}, and the
proof of proposition \ref{prop_a-un-and-quots} shows that the
sought identity follows from condition $\mathrm{(c)_{\bff}}$.
In particular, $\tau_A$ is a discrete subset of $A$, and
the natural map $\tau_A\to\tau_AA^\wedge$ is an isomorphism;
then, $\tau_A$ is (naturally identified with) a discrete
subset of $A^\wedge$ as well, and (i) implies that the natural
map $\tau_A\to\tau_{A^\wedge}$ is an isomorphism as well.

(iii)\ \
Consider now the special case where $J=bA$ for some $b\in A$.
Then we claim more precisely that $\Ann_{A^\wedge}(b^n)$ is the
topological closure of the image of $\Ann_A(b^n)$ in $A^\wedge$,
for every $n\in\N$. Indeed, let $x\in\Ann_{A^\wedge}(b^n)$; for
every $k\in\N$ we may find $x_k\in A$ and $y_k\in A^\wedge$
such that $x=x_k+b^ky_k$ in $A^\wedge$. Then
$b^nx_k\in(b^{n+k}A^\wedge)\cap A=b^{n+k}A$; say that
$b^nx_k=b^{n+k}z_k$ with $z_k\in A$. Hence the sequence
$(x_k-b^kz_k~|~k\in\N)$ lies in $\Ann_A(b^n)$ and converges
$b$-adically to $x$, whence the contention.
\end{remark}

\begin{corollary}\label{cor_pescato}
In the situation of \eqref{subsec_Urquart},
let $B$ be a ring, $\phi:B\to A$ a ring homomorphism,
$K\subset B$ and $J\subset I$ two ideals. Set
$M:=K\otimes_BA$, and denote by
$$
\psi:M\to A
$$
the $A$-linear map induced by $\phi$. Suppose moreover
that $A$ satisfies condition $\mathrm{(a)_\bff}$ of
\eqref{subsec_badabum}, and the image of $\psi$ is open
in the $I$-adic topology of $A$. Then we have :
\begin{enumerate}
\item
There exists $n\in\N$ such that $I^nM\cap\Ker\,\psi=0$.
\item
$A$ is $J$-adically complete and separated if and only
if the same holds for $M$.
\end{enumerate}
\end{corollary}
\begin{proof}(i): Notice first that $\Ker\,\psi=\Tor_1^B(B/K,A)$.
Especially, $\Ker\,\psi$ is naturally a $B/K$-module.
By assumption, there exists $k\in\N$ such that
$I^k\subset\Img\,\psi=\phi(K)\cdot A$; we deduce that
$\Ker\,\psi$ is also an $A/I^k$-module. Next, we consider
the short exact sequences
$$
\begin{aligned}
0\to &\, \Ker\,\psi\to M\to\Img\,\psi\to 0 \\
0\to &\, \Img\,\psi\to A\to A/\Img\,\psi\to 0.
\end{aligned}
$$
There follows exact sequences of inverse systems
$$
\begin{aligned}
T_\bullet:=&(\Tor^{A_r}_1(A_r/I^n_r,\Img\,\psi)~|~n\in\N)
\!\to\! K_\bullet:=(A/I^n\otimes_A\Ker\,\psi~|~n\in\N)\!\to\!
(M/\!I^nM~|~n\in\N) \\
T'_\bullet:=&(\Tor^{A_r}_2(A_r/I^n_r,A/\Img\,\psi)~|~n\in\N)
\to T_\bullet\to T''_\bullet:=(\Tor^{A_r}_1(A_r/I^n_r,A)~|~n\in\N).
\end{aligned}
$$
Now, since $A$ satisfies condition $\mathrm{(a)}_\bff$,
the inverse system $T''_\bullet$ is essentially zero,
and the same holds for $T'_\bullet$, by virtue of
lemma \ref{lem_ess-zero-Tors}; hence, $T_\bullet$ is
essentially zero, by lemma \ref{lem_inverse-Serre-subcat}.
Moreover, $K_n=\Ker\,\psi$ for every $n\geq k$, and the
transition maps $K_{n+1}\to K_n$ are the identities.
Then the assertion follows by a simple diagram chase : we
leave the details to the reader.

(ii): Since $\Img\,\psi$ is open in $A$ for the $I$-adic
topology, it is also open in $A$ for the $J$-adic topology,
and the latter agrees with the topology induced by the
$J$-adic topology of $A$; therefore, $A$ is $J$-adically
complete and separated if and only if the same holds for
$\Img\,\psi$. On the other hand, the $J$-adic topology on
$\Img\,\psi$ also agrees with the quotient topology induced
by the $J$-adic topology of $M$. Moreover, (ii) implies that
the $J$-adic topology of $M$ induces the discrete topology
on $\Ker\,\psi$. Then, proposition
\ref{prop_replaces-Mat-Th.8.1}(i,v) yields a short
exact sequence of $J$-adic completions :
$$
0\to\Ker\,\psi\to M^\wedge\to(\Img\,\psi)^\wedge\to 0
$$
which shows that $\Img\,\psi$ is $J$-adically complete
and separated if and only if the same holds for $M$,
whence (ii).
\end{proof}

\section{Complements of commutative algebra}
\label{chap_comm-algebra}
This chapter is a miscellanea of results of commutative
algebra that shall be needed in the rest of the treatise.

\subsection{Valuation theory}
This section is a selection of topics in valuation theory,
and will be complemented by the section \ref{sec_Spv-of-ring},
devoted to Huber's theory of the valuation spectrum.

\begin{definition}\label{def_ordered-group}
Recall that an {\em ordered abelian group} is a datum
$(\Gamma,\cdot,1,\leq)$ consisting of an abelian group
$(\Gamma,\cdot,1)$ and a total ordering $\leq$ on $\Gamma$
(see example \ref{ex_universe}(iii)) such that 
$$
\gamma\leq\gamma'
\quad\Rightarrow\quad
\gamma\cdot\gamma''\leq\gamma'\cdot\gamma''
\qquad
\text{for every $\gamma,\gamma',\gamma''\in\Gamma$}.
$$
\begin{enumerate}
\item
We denote by $\Gamma^+\subset\Gamma$ the submonoid of
all $\gamma\in\Gamma$ such that $\gamma\leq 1$.
\item
A subgroup $\Delta$ of the ordered abelian group $\Gamma$
is said to be {\em convex} if it satisfies the following
condition. If $\delta\in\Delta^+:=\Delta\cap\Gamma^+$, and
$\gamma\in\Gamma$ is any element with $1\geq\gamma\geq\delta$,
then $\gamma\in\Delta$. The {\em spectrum} of $\Gamma$ is
the set of all convex subgroups of $\Gamma$, denoted
$$
\Spec\,\Gamma.
$$
\item
The {\em convex rank} of $\Gamma$, denoted
$$
c.\rk(\Gamma)\in\N\cup\{\infty\}
$$
is the supremum over the lengths $r$ of the chains
$0\subset\Delta_1\subsetneq\cdots\subsetneq\Delta_r:=\Gamma$
of convex subgroups of $\Gamma$. The {\em rational rank}
of $\Gamma$ is 
$$
\rk(\Gamma):=\dim_\Q(\Gamma\otimes_\Z\Q)\in\N\cup\{\infty\}.
$$
\item
A morphism of ordered abelian groups is a group
homomorphism that is an order-preserving map for
the underlying totally ordered sets.
\end{enumerate}
\end{definition}

\begin{remark}\label{rem_ordered-gps}
Let $(\Gamma,\cdot,1,\leq)$ be any ordered abelian group.

(i)\ \
It is easily seen that $\Gamma$ is a torsion-free abelian
group. Moreover, we have quite generally
$$
c.\rk(\Gamma)\leq\rk(\Gamma)
$$
(details left to the reader).

(ii)\ \
Notice as well that there exists a unique structure of
ordered abelian group on $\Gamma\otimes_\Z\Q$ such that
the natural map $\Gamma\to\Gamma\otimes_\Z\Q$ is a morphism
of ordered groups. Namely, we have $\gamma\otimes q\leq 1$
in $\Gamma\otimes_\Z\Q$ whenever $\gamma\leq 1$ in $\Gamma$
and $q\geq 0$. Moreover, for every $q\in\Q$ we have
a well defined group endomorphism
$$
\Gamma\otimes_\Z\Q\to\Gamma\otimes_\Z\Q
\qquad
(\gamma\otimes r)\mapsto(\gamma\otimes r)^q:=\gamma\otimes qr.
$$

(iii)\ \
A subgroup $\Delta$ of $\Gamma$ is convex if and only
if there exists an ordered abelian group structure on
$\Gamma/\Delta$ such that the projection
$\Gamma\to\Gamma/\Delta$ is a morphism of ordered
abelian groups. Then the ordered abelian group structure
with this property is unique. Moreover, notice that the
intersection of an arbitrary family of convex subgroups
of $\Gamma$ is still convex. Especially, for every
subgroup $\Delta$ of $\Gamma$, there exists a minimal
convex subgroup of $\Gamma$ containing $\Delta$, called
the {\em convex hull} of $\Delta$; explicitly, it is the
subgroup of all elements $\gamma\in\Gamma$ such that there
exists $\delta\in\Delta^+$ with $\delta\leq\gamma\leq\delta^{-1}$.

(iv)\ \
There is a natural inclusion-reversing bijection
$$
\Spec\,\Gamma\isom\Spec\,\Gamma^+
\qquad
\Delta\mapsto\Gamma^+\setminus\Delta^+
$$
(notation of \eqref{subsec_sepc-of-monoid} and definition
\ref{def_ordered-group}(ii)). Indeed, it is easily seen
that $\Gamma^+\setminus\Delta^+$ is a prime ideal of
$\Gamma^+$, for every convex subgroup $\Delta$ of $\Gamma$.
Conversely, for any prime ideal $\fp\subset\Gamma^+$, the
subset $\Delta(\fp)^+:=\Gamma^+\setminus\fp$ is a submonoid,
and we claim that the subgroup $\Delta(\fp)$ generated
by $\Delta(\fp)^+$ is convex. Indeed, suppose that
$\gamma\in\Delta(\fp)^+$ and $\gamma'\in\Gamma^+$ is
any other element such that $\gamma\leq\gamma'$; then
$\beta:=\gamma\cdot\gamma'^{-1}\in\Gamma^+$, therefore
$\gamma'$ cannot lie in $\fp$, for otherwise the same
would hold for $\gamma=\beta\cdot\gamma'$. The contention
follows.

(v)\ \
Let $f:\Gamma\to\Gamma'$ be a morphism of ordered abelian
groups. Then $f$ induces a mapping
$$
\Spec\,\Gamma'\to\Spec\,\Gamma
\qquad
\Delta\mapsto f^{-1}\Delta.
$$
For instance, if $\Delta\subset\Gamma$ is any convex subgroup,
and $\pi:\Gamma\to\Gamma/\Delta$ the projection, then
$\Spec\,\pi:\Spec\,\Gamma/\Delta\to\Spec\,\Gamma$ is an
injective map, whose image is the subset of all convex
subgroups of $\Gamma$ containing $\Delta$. In the same vein,
notice that the inclusion map $i:\Gamma\to\Gamma\otimes_\Z\Q$
induces a bijection
$$
\Spec\,i:\Spec\,\Gamma\otimes_\Z\Q\isom\Spec\,\Gamma
$$
(details left to the reader). Moreover, if $f$ as above
is injective, then $\Spec\,f$ is surjective : indeed,
if $\Delta\subset\Gamma$ is any convex subgroup, let
$\Delta'$ be the convex hull of $f(\Delta)$ in $\Gamma'$;
then it is easily seen that $f^{-1}\Delta'=\Delta$.

(vi)\ \
If $c.\rk(\Gamma)=1$, we may find an injective morphism
of ordered abelian groups
$$
\rho:(\Gamma,\cdot,1,\leq)\to(\R,+,0,\leq).
$$
Indeed, pick an element $\gamma\in\Gamma$ with $\gamma>1$.
For every $\delta\in\Gamma$ and every positive integer $n$,
there exists a largest integer $k(n)$ such that
$\gamma^{k(n)}<\delta^n$. Then $(k(n)/n~|~n\in\N)$ is a
Cauchy sequence and we let
$\rho(\gamma):=\lim_{n\to\infty}k(n)/n$. One verifies easily
that $\rho$ is an order-preserving group homomorphism,
and since the convex rank of $\Gamma$ equals one, it
follows that $\rho$ is injective.

(vii)\ \
To any $\gamma\in\Gamma^+$ we may attach two convex subgroups
$$
o(\gamma):=\bigcap_{n\in\N}
\{\delta\in\Gamma~|~\gamma\leq\delta^n\leq\gamma^{-1}\}
\qquad\text{and}\qquad
O(\gamma):=\bigcup_{n\in\N}
\{\delta\in\Gamma~|~\gamma^n\leq\delta\leq\gamma^{-n}\}.
$$
Clearly $o(\gamma)\subset O(\gamma)$, and if
$\gamma\neq 1$, the quotient $Q(\gamma):=O(\gamma)/o(\gamma)$
has convex rank one (for its natural ordered group
structure as in (iii)). Indeed, suppose that
$\delta,\mu$ are two elements of $O(\gamma)$ with
classes $\bar\delta,\bar\mu\in Q(\gamma)_+$ and
$\bar\delta\neq 1$; the assumptions mean that there
exist $n,k\in\N$ such that $\delta^k<\gamma$ and
$\gamma^n\leq\mu\leq\gamma^{-n}$. Then
$\delta^{kn}\leq\mu$; thus, the only convex subgroups
of $Q(\gamma)$ are $Q(\gamma)$ and $\{1\}$, as claimed.

(viii)\ \
Standard examples of ordered abelian groups are
$(\Q,+,0,\leq)$ and  $(\R,+,0,\leq)$; the
latter is also isomorphic to the ordered group
$(\R_{>0},\cdot,1,\leq)$ of strictly positive
real numbers. We shall also use the notation :
$\Q_{>0}:=\Q\cap\R_{>0}$, $\R_+:=\R_{>0}\cup\{0\}$ and
$\Q_+:=\Q_{>0}\cup\{0\}$.
\end{remark}

\begin{definition}\label{def_valuation}
Let $A$ be any ring, $(\Gamma,\cdot,1,\leq)$ any ordered
abelian group; we extend the ordering and the composition
law of $\Gamma$ to $\Gamma_\circ:=\Gamma\cup\{0\}$
(cp. remark \ref{rem_apparent-reasons}(i)), by
the rule
$$
0<\gamma
\qquad\text{and}\qquad
0\cdot\gamma=0=\gamma\cdot 0
\qquad
\text{for every $\gamma\in\Gamma$}.
$$
Moreover, for every $q\in\Q_+$ we extend to
$(\Gamma\otimes_\Z\Q)_\circ$ the $q$-th power operation
of remark \ref{rem_ordered-gps}(ii), by declaring that
$0^0:=1$ and $0^q:=0$ for every strictly positive $q\in\Q$.

(i)\ \
A {\em $\Gamma$-valued semi-norm} on $A$ is a mapping
$$
v:A\to\Gamma_\circ
$$
such that :
\begin{itemize}
\item
$v(0)=0$ and $v(1)\leq 1$.
\item
$v(a-b)\leq\max(v(a),v(b))$ for every $a,b\in A$.
\item
$v(ab)\leq v(a)\cdot v(b)$ for every $a,b\in A$.
\end{itemize}
The {\em value group} of $v$ is the subgroup
$$
\Gamma_{\!v}
$$
of $\Gamma$ generated by $\Img\,(v)\setminus\{0\}$.
The {\em rank} of $v$ is the convex rank of its
value group.

(ii)\ \
We say that the semi-norm $v$ is {\em power-multiplicative}
if we have :
$$
v(a^n)=v(a)^n
\qquad
\text{for every $a\in A$ and every integer $n>0$}.
$$

(iii)\ \
We say that the semi-norm $v$ is a {\em valuation} if we have :
$$
v(1)=1
\qquad\text{and}\qquad
v(ab)=v(a)\cdot v(b)
\qquad\text{for every $a,b\in A$}.
$$

(iv)\ \
We say that $v$ is a {\em norm} if $v(a)\neq 0$
for every $a\neq 0$ in $A$.

(v)\ \
We say that the semi-norm $v$ is {\em real-valued}
if its value group is a subgroup of $\R_{>0}$.

(vi)\ \
Let $\Gamma$ and $\Gamma'$ be two ordered abelian groups,
$v$ (resp. $v'$) a $\Gamma$-valued (resp. $\Gamma'$-valued)
semi-norm on $A$. We say that $v$ is {\em equivalent} to
$v'$ if there exist an ordered abelian group $\Gamma''$
and a semi-norm $v'':A\to\Gamma''$, together with injective
maps of ordered abelian groups $\phi:\Gamma''\to\Gamma$ and
$\phi':\Gamma''\to\Gamma'$ such that
$\phi_\circ\circ v''=v$ and $\phi'_\circ\circ v''=v'$
(notation of remark \ref{rem_apparent-reasons}(i)).
\end{definition}

\begin{remark}\label{rem_semi-norm}
Let $A$ and $\Gamma$ be as in definition \ref{def_valuation},
and $|\cdot|_A$ any $\Gamma$-valued semi-norm on $A$.

(i)\ \
Notice that $|1|_A=|1\cdot 1|_A\leq|1|^2_A\leq|1|_A$, whence
$|1|_A=|1|^2_A$, so that either $|1|_A=1$ or $|1|_A=0$.
In case $|1|_A=0$, then $|a|_A=|a\cdot 1|_A\leq|a|_A\cdot 0=0$,
and in this case $|\cdot|_A$ is the trivial semi-norm with
constant value $0$.

(ii)\ \
We have $|-a|_A=|a|_A$ for every $a\in A$. Indeed,
$|-a|_A=|0-a|_A\leq\max(|0|_A,|a|_A)=|a|_A$; after replacing
$a$ by $-a$ we get also $|a|_A\leq|-a|_A$, whence the claim.

(iii)\ \
Let $a,b\in A$ be any two elements such that $|a|_A>|b|_A$.
Then we have $|a+b|_A=|a|_A$. Indeed, suppose that
$|a+b|_A<|a|_A$; taking into account (ii) we get
$$
|a|_A\leq\max(|a+b|_A,|-b|_A)<|a|_A
$$
which is absurd.

(iv)\ \
Let $S\subset A$ be any multiplicative subset such
that $|s|_A\neq 0$ for every $s\in S$. If $|\cdot|_A$
is a valuation, there exists a unique valuation
$$
|\cdot|_{S^{-1}A}:S^{-1}A\to\Gamma_\circ
$$
whose composition with the localization map $A\to S^{-1}A$
agrees with $|\cdot|_A$.

(v)\ \
Notice that the subsets
$$
\Ker\,(|\cdot|_A):=\{a\in A~|~|a|_A=0\}\subset
(A,|\cdot|_A)^+:=\{a\in A~|~|a|_A\leq 1\}
$$
are respectively an ideal and a subring of $A$.
Clearly $|\cdot|_A$ is the composition of the projection
$$
A\to A/\Ker\,(|\cdot|_A)
$$
and a unique $\Gamma$-valued norm on $A/\Ker\,(|\cdot|_A)$.
Moreover, if $|\cdot|_A$ is a valuation, then
$\Ker\,(|\cdot|_A)$ is a prime ideal of $A$ called the
{\em support}\/ of $|\cdot|_A$, and we define the
{\em residue field\/} of $|\cdot|_A$ as
$$
\kappa(|\cdot|_A):=\Frac\,A/\Ker\,(|\cdot|_A).
$$
In light of (iv), in this case $|\cdot|_A$ factors
uniquely through the projection $A\to\kappa(|\cdot|_A)$
and a valuation $|\cdot|_\kappa$ on $\kappa(|\cdot|_A)$,
called the {\em residual valuation of\/ $|\cdot|_A$}.
Conversely, if $\fp$ is any prime ideal of $A$, and
$|\cdot|_{\kappa(\fp)}$ is any valuation on the residue
field $\kappa(\fp)$ of $\fp$, then the composition of
$|\cdot|_{\kappa(\fp)}$ with the projection $A\to\kappa(\fp)$
is a valuation on $A$.

(vi)\ \
With the notation of (v), we claim that if $|\cdot|_A$ is
a power-multiplicative semi-norm, then the subring
$(A,|\cdot|_A)^+$ is integrally closed in $A$. Indeed,
let $a\in A$ be any element that satisfies an identity
of the type
$$
a^n+b_1a^{n-1}+\cdots+b_n=0
\qquad
\text{with $b_1,\dots,b_n\in(A,|\cdot|_A)^+$}.
$$
It follows that $|a|^n_A=|a^n|_A\leq
\max(|b_i|_A\cdot|a|_A^{n-i}~|~i=1,\dots,n)$.
Hence, say that $|a|^n_A\leq|b_i|\cdot|a|^{n-i}_A$ for
some $i\leq n$; if $|a|=0$, obviously $a\in(A,|\cdot|_A)^+$.
Otherwise, we deduce $|a|^i_A\leq|b_i|_A\leq 1$, whence
$|a|_A\leq 1$, and the claim follows.

(vii)\ \
Let $\fp\subset A$ be any prime ideal. Then there exists
a unique (up to equivalence) rank zero valuation on $A$
with support $\fp$, defined by the rule :
$$
|x|:=\left\{\begin{array}{ll}
       0 & \text{if $x\in\fp$} \\
       1 & \text{otherwise}.
       \end{array}\right.
$$
Any valuation of this type is called a {\em trivial valuation}
on $A$.
\end{remark}

\begin{example}\label{ex_toric-valuations}
(i)\ \
Let $A$ be any ring, $P$ a monoid, $(\Gamma,\cdot,1,\leq)$
an ordered abelian group, and $\phi:P\to\Gamma_\circ$ a
morphism of monoids. We deduce a mapping (notation of
\eqref{subsec_mon-to-algs})
$$
v_\phi:A[P]\to\Gamma_\circ
\qquad
\sum_{x\in P}a_x\cdot x\mapsto\max(\phi(x)~|~x\in P,\ a_x\neq 0)
$$
(where the maximum of the empty set is taken to be
$0\in\Gamma_\circ$, so that $w_\phi(0)=0$). It is easily seen
that $v_\phi$ is a $\Gamma$-valued seminorm on $A[P]$. Also,
$v_\phi$ is a norm if and only if $\phi(P)\subset\Gamma$.

(ii)\ \
Suppose moreover that $A$ is a reduced ring, $P$ is integral
and $P^\gp$ is a torsion-free abelian group; then $v_\phi$ is
valuation. Indeed, consider the first the case where $A$ is
a domain; then the same holds for $A[P]$, since the latter
is a subring of $A[P^\gp]$, which in turn is the filtered
union of its subrings $A[G]$, with $G$ ranging over all
free abelian subrings of finite rank contained in $P^\gp$.
Now, let
$$
a:=\sum_{x\in P}a_x\cdot x
\qquad
b:=\sum_{x\in P}b_x\cdot x
$$
be any two non-zero elements of $A[P]$, and set $\alpha:=v_\phi(a)$
and $\beta:=v_\phi$; let also $a':=\sum_{\phi(x)=\alpha}a_x\cdot x$,
and $b':=\sum_{\phi(x)=\beta}b_x\cdot x$. It is easily seen that
$v_\phi(ab)=v_\phi(a'b')$; on the other hand we have $a'b'\neq 0$
since $A[P]$ is a domain, so we must have $v_\phi(a'b')=\alpha\beta$,
whence the assertion.
If $A$ is any reduced ring, we have a natural inclusion
$A\to\prod_\fp A_\fp$, where $\fp$ ranges over all minimal
prime ideals of $A$, and each factor $A_\fp$ is a field.
Then we get likewise an inclusion $A[P]\to\prod_\fp A_\fp[P]$,
and the assertion is easily reduced to the foregoing case :
details left to the reader.
\end{example}

\sset\subsubsection{}\label{subsec_uniform-semi-normed}
Real-valued semi-norms, of course, occur frequently in
algebraic and analytic questions. The datum $(A,|\cdot|_A)$
of a ring and a real valued semi-norm $|\cdot|_A$ is usually
called a {\em (real valued) semi-normed ring}. If $|\cdot|_A$
is power-multiplicative, one also says that $(A,|\cdot|_A)$
is a {\em uniform semi-normed ring}. A {\em morphism of
semi-normed rings} $f:(A,|\cdot|_A)\to(B,|\cdot|_B)$ is a
ring homomorphism $f:A\to B$ for which there exists a real
number $C\geq 0$ such that
\set\begin{equation}\label{eq_morph-of-semi-normed}
|f(a)|_B\leq C\cdot|a|_A
\qquad
\text{for every $a\in A$}.
\end{equation}
Clearly a composition of morphisms of semi-normed ring
is again a morphism of semi-normed rings, hence the
semi-normed rings and their morphisms form a category
$\cS$, and we let $\su.\cS\subset\cS$ be the full
subcategory whose objects are the uniform semi-normed rings.
Notice that if $f$ verifies \eqref{eq_morph-of-semi-normed}
and $(B,|\cdot|_B)$ is a uniform semi-normed ring, then
we have
$$
|f(a)|_B^n=|f(a^n)|_B\leq C\cdot|a^n|\leq C\cdot|a|^n
\qquad
\text{for every $n\in\N$ and $a\in A$}
$$
whence $|f(a)|_B\leq C^{1/n}\cdot|a|_A$ for every such
$n$ and $a$; especially, \eqref{eq_morph-of-semi-normed}
holds with $C\leq 1$.

\begin{lemma}\label{lem_normalized-Samuel}
Let $(A,|\cdot|)$ be any (real valued) semi-normed ring.
Then the following holds :
\begin{enumerate}
\item
The sequence $(|a^n|^{1/n}~|~n\in\N)$ converges in
$\R$ for every $a\in A$. We set
$$
|a|^*:=\lim_{n\to+\infty}|a^n|^{1/n}.
$$ 
\item
The mapping $|\cdot|^*:A\to\R$ is a power-multiplicative
semi-norm on $A$.
\item
$|a|^*\leq|a|$ for every $a\in A$.
\item
The rule : $(A,|\cdot|)\mapsto(A,|\cdot|^*)$ defines
a left adjoint for the inclusion functor $\su.\cS\to\cS$.
\end{enumerate}
\end{lemma}
\begin{proof}(i): Fix $a\in A$, and set $\nu_n:=|a^n|^{1/n}$
for every integer $n>0$. Suppose first that $\nu_k=0$
for some $k\in\N$; it follows easily that $\nu_n=0$ for
every $n\geq k$, whence the assertion, in this case.
Hence, we may assume that $\nu_n>0$ for every $n\in\Z$;
we notice that
$$
\nu_{kn}\leq\nu_n
\qquad
\text{for every $k,n>0$}
$$
and therefore
$$
\nu_{kn+r}\leq(|a^{kn}|\cdot|a^r|)^{1/(kn+r)}\leq
\nu_n^{1-1/(k+1)}\cdot|a^r|^{1/(nk+r)}
$$
for every $k,n,r\in\N$ with $0\leq r<n$ and $k>0$.
Now, for a given integer $n>0$ and a real number $\eps>0$,
pick $k\in\N$ large enough, so that $|a^r|^{1/(nk)}>1-\eps$
for every $r=0,\dots,n-1$, and $1/(k+1)<\eps$; it follows
that
$$
\nu_t\leq\nu_n^{1-\eps}\cdot(1-\eps)
\qquad
\text{for every integer $t\geq kn$}.
$$
Consequently
$$
\limsup_{n\to+\infty}\nu_n=\liminf_{n\to+\infty}\nu_n
$$
whence the assertion.

(ii): Fix $a,b\in A$ and a real number $\eps>0$;
by (i), there exists $k_0\in\N$ such that
$$
|a^k|^{1/k}<|a|^*\cdot(1+\eps)
\qquad\text{and}\qquad
|b^k|^{1/k}<|b|^*\cdot(1+\eps)
\qquad
\text{for every $k>k_0$}.
$$
Pick an integer $t\geq 2k_0$, and large enough, so that
$$
|a^i|^{1/n},|b^i|^{1/n}<1+\eps
\qquad
\text{for every $n\geq t$ and every $i\leq k_0$}.
$$
Now, we have :
$$
|(a+b)^n|^{1/n}\leq\max
\{|a^i|^{1/n}\cdot|b^{n-i}|^{1/n}~|~i=0,\dots,n\}.
$$
On the other hand, we may assume that $|a|^*\geq|b|^*$;
set also
$$
M:=\max((|a|^*\cdot(1+\eps))^{1-\eps},|a|^*\cdot(1+\eps)).
$$
Then, for every $n\in\N$ such that $n\geq t$ and
$k_0/n\leq\eps$, we get
$$
|a^i|^{1/n}\cdot|b^{n-i}|^{1/n}\leq
\left\{\begin{array}{ll}
          |a|^*\cdot(1+\eps) & \text{whenever $i,n-i>k_0$} \\
           M\cdot(1+\eps)    & \text{otherwise}.
       \end{array}\right.
$$
We immediately deduce that $|a+b|^*\leq\max(|a|^*,|b|^*)$
for every $a,b\in A$.

The inequality $|ab|^*\leq|a|^*\cdot|b|^*$ for every $a,b\in A$
is an immediate consequence of the corresponding property for
$|\cdot|$, and the power-multiplicative condition follows
easily from the definition of $|\cdot|^*$.

(iii) follows easily by remarking that $|a^n|\leq|a|^n$ for
every $a\in A$.

(iv): Let $(B,|\cdot|_B)$ be a uniform semi-normed ring, and
$f:(A,|\cdot|)\to(B,|\cdot|_B)$ a morphism of semi-normed rings,
so that \eqref{eq_morph-of-semi-normed} holds for some $C\geq 0$.
We get : $|f(a)|_B=|f(a^n)|^{1/n}_B\leq C^{1/n}\cdot|a^n|^{1/n}$
for every $a\in A$ and $n\in\N$. After taking the limit for
$n\to+\infty$, it follows that $|f(a)|_B\leq|a|^*$ for every
$a\in A$, so that $f:(A,|\cdot|^*)\to(B,|\cdot|_B)$ is a morphism
in $\su.\cS$. The assertion follows immediately.
\end{proof}

\begin{example}\label{ex_Samuel}
(i)\ \
Let $A$ be any ring, $A_0\subset A$ a subring, $I\subset A_0$
an ideal. We define, for every $n\in\N$, the $A_0$-submodule
of $A$
$$
I^{-n}:=\{a\in A~|~aI^n\subset A_0\}
$$
and we consider the {\em order function associated with $I$}
$$
\nu_I:A\to\Z\cup\{\pm\infty\}
\qquad
a\mapsto\sup\{n\in\Z~|~a\in I^n\}
$$
(where the supremum of the empty set is taken to be $-\infty$).
It is easily seen that
\begin{itemize}
\item
$\nu_I(a+b)\geq\min(\nu_I(a),\nu_I(b))$ for every $a,b\in A$
\item
$\nu_I(ab)\geq\nu_I(a)+\nu_I(b)$ for every $a,b\in A$
\item
$\nu_I(a)\geq 0$ if and only if $a\in A_0$.
\end{itemize}
Let also
$$
A(I):=\{a\in A~|~\nu_I(a)\neq-\infty\}.
$$
It is easily seen that $A(I)$ is a subring of $A$. Moreover,
fix any real number $\rho\in]0,1[$, and set
$$
|a|_I:=\rho^{\nu_I(a)}\in\R_{\geq 0}.
\qquad
\text{for every $a\in A(I)$}
$$
(with the convention that $\rho^{+\infty}:=0$). It is clear
that $|\cdot|_I$ is a semi-norm on $A(I)$.
The corresponding power-multiplicative semi-norm
$|\cdot|^*_I$ provided by lemma \ref{lem_normalized-Samuel}
is sometimes called the {\em asymptotic Samuel function}
of $I$.

(ii)\ \
In the situation of (i), let $J\subset A_0$ be another
ideal, with $I\subset J$. Then it is easily seen that,
for every $a\in A$, we have either $\nu_J(a)\geq\nu_I(a)\geq 0$,
or else $0>\nu_I(a)\geq\nu_J(a)$. Therefore
$$
A(J)\subset A(I)
$$
and the two foregoing conditions can be unified in
the following assertion. For every $a\in A(J)$ we have
either $\nu_J(a)=0$, in which case $\nu_I(a)=0$ as well,
or else
$$
0\leq\nu_I(a)/\nu_J(a)\leq 1.
$$
We deduce a corresponding inequality for asymptotic Samuel
functions : namely, for every $a\in A(J)$ with $|a|^*_J\neq 0$
we have either $|a|^*_J=1$, in which case $|a|^*_I=1$ as well,
or else
$$
0\leq\frac{\log|a|^*_I}{\log|a|^*_J}\leq 1.
$$

(iii)\ \
For instance, for any integer $n>0$ we have $A(I)=A(I^n)$ and
$$
\nu_{I^n}(a^n)\geq\nu_I(a)\geq n\cdot\nu_{I^n}(a)
\qquad
\text{for every $a\in A$}
$$
from which it follows easily that
$$
|a|^*_{I^n}=(|a|^*_I)^{1/n}
\qquad
\text{for every $a\in A(I)$}.
$$

(iv)\ \
Let $I$ and $I'$ two ideals of $A_0$, and suppose
that the $I$-adic topology on $A_0$ agrees with the
$I'$-adic topology. Combining (ii) and (iii) we deduce
that
$$
A(I)=A(I')
$$
and for every $a\in A(I)$ we have $|a|^*_I=1$ (resp.
$|a|^*_I=0$) if and only if $|a|_{I'}^*=1$ (resp.
$|a|^*_{I'}=0$); also, there exists a real number
$C\geq 0$ (independent of $a$) such that
$$
0\leq\frac{\log|a|^*_I}{\log|a|^*_{I'}},
\frac{\log|a|^*_{I'}}{\log|a|^*_I}\leq C
\qquad
\text{whenever $|a|^*_I\neq 1,0$}.
$$

(v)\ \
We notice as well that the asymptotic Samuel function
is independent of the subring $A_0$. Namely, in the
situation of (i), suppose that $A'_0$ is another subring
of $A$ such that $I$ is also an ideal of $A'_0$. Then
we may set
$$
I'^n:=I^n
\qquad
\text{for every $n>0$, and}
\qquad
I'^{-k}:=\{a\in A~|~aI^k\subset A'_0\}
\qquad
\text{for every $k\geq 0$}
$$
(where $I^0:=A_0$). Using the system of powers $(I^n~|~n\in\Z)$
we define the function $\nu_I$ and the subring $A(I)$ as in (i),
and with the system of powers $(I'^n~|~n\in\Z)$ we may define
likewise the function $\nu_{I'}$ and the subring $A(I')$. After
fixing $\rho\in]0,1[$, we then may define correspondingly the
semi-norms $|\cdot|_I$ and $|\cdot|_{I'}$, as well as the
associated power multiplicative semi-norms $|\cdot|^*_I$ and
$|\cdot|^*_{I'}$. Obviously, for every $a\in A$ we have
$\nu_I(a)>0$ if and only $\nu_{I'}(a)>0$, and
$\nu_I(a)=\nu_{I'}(a)$ for every $a\in A$ such that $\nu_I(a)>0$.
Moreover, suppose that $\nu_I(a)=-n$ for some $n\in\N$; this
means especially that $aI^n\subset A_0$, so that
$aI^{n+1}\subset I\subset A'_0$, whence $\nu_I(a)\geq-n-1$.
It follows easily that
$$
A(I)=A(I')
\qquad\text{and}\qquad
|a|^*_I=|a|^*_{I'}
\qquad
\text{for every $a\in A(I)$}.
$$

(vi)\ \
It follows especially from (v) that the subring
$(A(I),|\cdot|^*_I)^+$ of $A(I)$ is independent of $A_0$
(notation of remark \ref{rem_semi-norm}(v)), and
moreover, if $B\subset A$ is any subring such that
$I$ is also an ideal of $B$, then $B\subset(A(I),|\cdot|^*_I)^+$.
More generally, we have :
\end{example}

\begin{proposition} In the situation of example
{\em\ref{ex_Samuel}(i)}, consider any subring $B$ of $A$ with
$$
I\subset B\subset(A(I),|\cdot|^*_I)^+.
$$
Then the following holds :
\begin{enumerate}
\item
$A(I)\subset A(IB)$\ \ and\ \ $|a|^*_I=|a|^*_{IB}$\ \
for every $a\in(A(I),|\cdot|^*_I)^+$.
\item
Suppose moreover that $I$ is a finitely generated ideal of
$A_0$. Then $A(I)=A(IB)$\ \ and\ \ $|a|^*_I=|a|^*_{IB}$\ \
for every $a\in A(I)$.
\end{enumerate}
\end{proposition}
\begin{proof}(i): Arguing as in example \ref{ex_Samuel}(v),
we easily see that $\nu_{IB}(a)\geq\nu_I(a)$ for every
$a\in A$ such that $\nu_I(a)>0$, and
$\nu_{IB}(a)\geq\nu_I(a)-1$ if $\nu(a)\leq 0$. This already
implies the stated inclusion of subrings, as well as the
inequality
\set\begin{equation}\label{eq_vegan}
|a|^*_{IB}\leq|a|^*_I
\qquad
\text{for every $a\in A(I)$}.
\end{equation}
Now, suppose that $|a|^*_{IB}=\rho^c$ for some real number
$c>0$. Set $n(k):=\nu_{IB}(a^k)$ for every $k\in\N$. It
follows that $\lim_{k\to+\infty}n(k)/k=c$, and especially
$n(i)>0$ for every sufficiently large $i\in\N$. It
suffices to check that $|a^i|^*_I\leq\rho^{ic}$, so
we may replace $a$ by $a^i$ for such sufficiently large $i$,
and assume from start that $\nu_{IB}(a)\geq 0$. In this
situation, for every $k\in\N$ we may then write
$a^k=\sum_{j=1}^ra_ib_i$, for some $a_1,\dots,a_r\in I^{n(k)}$
and $b_1,\dots,b_r\in B$. Since by assumption $|b_i|^*_I\leq 1$
for every $i\leq r$, taking into account lemma
\ref{lem_normalized-Samuel}(ii,iii) we deduce that
$$
|a^k|^*_I\leq\sup(|x|^*_I~|~x\in I^{n(k)})\leq
\sup(|x|_I~|~x\in I^{n(k)})\leq\rho^{n(k)}
$$
whence $|a|^*_I\leq\rho^{n(k)/k}$ for every $k\in\N$,
and therefore $|a|^*_I\leq|a|^*_{IB}$.

(ii): Indeed, we know already that $A(I)\subset A(IB)$,
and to show the converse inclusion, pick any finite system
of generators $(t_1,\dots,t_r)$ for $I$, and let
$a\in A(IB)$ be any element. This means that
$a\cdot(IB)^n\subset B\subset A(I)$ for some $n\in\N$.
Hence we may find $m\in\N$ such that
$$
a\cdot t^{k_1}_1\cdots t^{k_r}_r\cdot I^m\subset A_0
\qquad
\text{for every $(k_1,\dots,k_r)\in\N^r$ with
$k_1+\cdots+k_r=n$}
$$
whence $a\cdot I^{m+n}\subset A_0$, so $a\in A(I)$, as
stated. Next, to conclude, in light of \eqref{eq_vegan},
it suffices to show that $|a|^*_I\leq|a|^*_{IB}$, for
every $a\in A(I)$. However, in case $|a|^*_{IB}=\rho^c$
for some real number $c>0$, we may argue as in the proof
of (i) to deduce the sought inequality. We may therefore
assume that $|a|^*_{IB}\leq\rho^{-c}$ for some real number
$c>0$, and we have to show that $|a|^*_I\leq\rho^{-c}$ as
well. Define $n(k)$ as in the proof of (i), for every
$k\in\N$. Notice that, if $n(k)>0$ for some $k\in\N$, we
have $|a^k|_{IB}\leq\rho$, whence $|a|^*_{IB}\leq\rho^{1/k}$,
in which case the assertion is already known. Hence, we
may further assume that $n(k)\leq 0$ for every $k\in\N$,
in which case
$$
\lim_{k\to+\infty}-n(k)/k\leq c
\qquad\text{and}\qquad
a^k\cdot(IB)^{-n(k)}\subset B
\qquad
\text{for every $k\in\N$}.
$$
Pick any integers $p,q>0$ such that $c<p/q$; it
follows easily that
$$
a^{nq}\cdot(IB)^{np-1}\subset B
\qquad
\text{for every sufficiently large $n\in\N$}.
$$
Fix such a sufficiently large $n\in\N$; we deduce
$$
x(k_\bullet):=a^{nq}\cdot t_1^{k_1}\cdots t_r^{k_r}\in IB
\qquad
\text{for every
$k_\bullet:=(k_1,\dots,k_r)\in\N^r$ with $k_1+\cdots+k_r=np$}
$$
in which case $|x(k_\bullet)|^*_{IB}<1$, and therefore
$|x(k_\bullet)|^*_I<1$ as well, by the foregoing, for
every such $k_\bullet$. Then we may find $m\in\N$
such that
$$
a^{nmq}\cdot t_1^{mk_1}\cdots t_r^{mk_r}\in A_0
\qquad
\text{for every
$k_\bullet:=(k_1,\dots,k_r)\in\N^r$ with $k_1+\cdots+k_r=np$}.
$$
However, let $J\subset A_0$ be the ideal generated
by $t_1^m,\dots,t_r^m$, and notice that
$$
I^{(m-1)(r+nps)+1}\subset J^{nps}
\qquad
\text{for every $s\in\N$}.
$$
We conclude that
$$
a^{nmqs}\cdot I^{(m-1)(r+nps)+1}\subset A_0
\qquad
\text{for every $s\in\N$}
$$
from which we see that
$$
\lim_{s\to+\infty}-\nu_I(a^{nmqs})/(nmqs)\leq
\lim_{s\to+\infty}((m-1)(r+nps)+1)/(nmqs)<p/q
$$
and since the rational number $p/q$ can be taken arbitrary
close to $c$, the sought inequality follows.
\end{proof}

\begin{definition}\label{def_dominant}
Let $F$ be a field.
\begin{enumerate}
\item
We say that a subring $V\subset F$ is a {\em valuation ring
of $F$} if the following holds. For every $x\in F^\times$
we have either $x\in V$ or $x^{-1}\in V$.
\item
Let $(A,\fm_A)$, $(B,\fm_B)$ be two local subrings of $F$.
We say that $B$ {\em dominates} $A$ if $A\subset B$ and
$A\cap\fm_B=\fm_A$.
\end{enumerate}
\end{definition}

\begin{remark}\label{rem_valuations}
(i)\ \
Let $K$ be any field; notice that any semi-norm, and
especially, any valuation on $K$ is necessarily a norm.

(ii)\ \
Let $|\cdot|_K$ be any such valuation, and let
$K^+:=(K,|\cdot|_K)^+$, as in remark \ref{rem_semi-norm}(v).
Then it is easily seen that $K^+$ is a valuation ring of
$K$. Also, the invertible elements of of $K^+$ are precisely
those $x\in K$ such that $|x|_K=1$, and the set
$\{x\in K~|~|x|_K<1\}$ is the unique maximal ideal of
$K^+$; especially, $K^+$ is a local domain. 

(iii)\ \
Furthermore, $K^+$ is integrally closed. Indeed, let
$x\in K$ be integral over $K^+$, so that
$x^n+y_1x^{n-1}+\cdots+y_n=0$ for some $y_1,\dots,y_n\in K^+$,
and suppose, by way of contradiction, that $x\notin K^+$,
{\em i.e.} that $|x|_K>1$; it follows easily that
$|x^n|_K>|y_ix^{n-i}|_K$ for every $i=1,\dots,n$,
so that remark \ref{rem_semi-norm}(iii) yields
$|x^n+y_1x^{n-1}+\cdots+y_n|_K=|x^n|_K\neq 0$, which
is absurd.

(iv)\ \
Conversely, if $V$ is any valuation ring of $K$, set
$\Gamma:=K^\times/V^\times$. Thus, $\Gamma$ is an abelian
group, and we define an ordering on $\Gamma$ as follows.
For given classes $\bar x,\bar y\in\Gamma$, we declare
that $\bar x\leq\bar y$ if and only if $y^{-1}x\in V$.
Then it is easily seen that $(\Gamma,\leq)$ is an ordered
abelian group, and the projection $K^\times\to\Gamma$
is the restriction to $K^\times$ of a well defined
valuation $|\cdot|_K$ of $K$, such that $V=(K,|\cdot|_K)^+$,
so the maximal ideal $\fm_V$ of $V$ and the subset $V^\times$
are also described as in (ii) (details left to the reader).

(v)\ \
If the valuation $|\cdot|_K$ is a surjective map,
we say that $(K,|\cdot|_K)$ is a {\em valued field}.
In view of (iv), we see that a valuation ring of $K$
is the same as the datum of the equivalence class (in
the sense of definition \ref{def_valuation}(vi)) of a
valued field $(K,|\cdot|_K)$.

(vi)\ \
Let $A$ be any ring, $|\cdot|$ and $|\cdot|'$ two valuations
on $A$; in light of (v), we see that the following conditions
are equivalent :
\begin{enumerate}
\alphaenu
\item
the valuations $|\cdot|$ and $|\cdot|'$ are equivalent
\item
$|\cdot|$ and $|\cdot|'$ have the same
support and $(\kappa(|\cdot|),|\cdot|_\kappa)^+=
(\kappa(|\cdot|'),|\cdot|'_\kappa)^+$
\item
for every $a,b\in A$ we have $|a|\leq|b|$ if and only
if $|a|'\leq|b|'$.
\end{enumerate}

(vii)\ \
Let $(K,|\cdot|)$ be any valued field, and $\Gamma$ the
value group of $|\cdot|$. There is an inclusion-reversing
bijection
$$
\Spec\,\Gamma\isom\Spec\,K^+
$$
that assigns to every convex subgroup $\Delta\subset\Gamma$
the prime ideal
$$
\fp(\Delta):=\{x\in K^+~|~\gamma>|x| \ 
            \text{for every $\gamma\in\Delta$}\}
$$
and whose inverse assigns, to every prime ideal $\fp\subset K^+$,
the convex subgroup
$$
\Delta(\fp):=\{\gamma\in\Gamma~|~\gamma,\gamma^{-1}>|x| \ 
                         \text{for every $x\in\fp$}\}.
$$
Indeed, it is easily seen that $\fp(\Delta)$ is a prime ideal,
and we check that $\Delta(\fp)$ is a convex subgroup. First,
say that $\gamma,\gamma'\in\Delta(\fp)$; if either $\gamma>1$
or $\gamma'>1$, it is clear that $\gamma\cdot\gamma'>|x|$
for every $x\in\fp$. If $\gamma,\gamma'\leq 1$, we may find
$y,y'\in K^+$ such that $\gamma=|y|$ and $\gamma'=|y'|$,
and by assumption, neither $y$ nor $y'$ lie in $\fp$, so also
$y\cdot y'\notin\fp$, and therefore $\gamma\cdot\gamma'>|x|_K$
for every $x\in\fp$. Hence, $\Delta(\fp)$ is a subgroup of
$\Gamma$, and the convexity of $\Delta(\fp)$ is obvious. Next,
set $\fq:=\fp(\Delta)$; directly from the definitions we see
that $\Delta\subset\Delta(\fq)$.
Conversely, let $\gamma\in\Gamma^+\setminus\Delta$; then
$\gamma=|y|$ for some $y\in K^+$, and since $\gamma<\delta$
for every $\delta\in\Delta$, we get $y\in\fq$, so that
$\gamma\notin\Delta_\fq$. Hence, $\Delta=\Delta(\fq)$, and
especially, the rule $\fp\mapsto\Delta(\fp)$ is surjective
from $\Spec\,K^+$ to $\Spec\,\Gamma$; but it is easily
seen that this rule is also injective, whence the contention.
\end{remark}

\begin{definition}\label{def_valuation-toplog}
Let $(K,|\cdot|_K)$ be any valued field, and $\Gamma$
the value group of $|\cdot|_K$.
\begin{enumerate}
\item
For every $\gamma\in\Gamma^+$, set
$U_\gamma:=\{x\in K^+~|~|x|_K<\gamma\}$.
The {\em valuation topology} of $K^+$ is the unique
linear topology $\cT_{K^+}$ that admits the family
$\{U_\gamma~|~\gamma\in\Gamma^+\}$ as a fundamental
system of open neighborhoods of $0$.
\item
In view of remark \ref{rem_I_bullet-adic}(i) there exists
a unique ring topology $\cT_K$ on $K$ such that :
\begin{enumerate}
\item
$K^+$ is an open subring of $K$
\item
$\cT_K$ induces on $K^+$ the valuation topology $\cT_{K^+}$.
\end{enumerate}
We call $\cT_K$ the {\em valuation topology} of $K$.
\item
We say that $(K,|\cdot|_K)$ is a {\em Tate valued field},
if the valuation topology of $K$ is f-adic and not discrete.
\item
Let $A$ be any ring, and $v$ any valuation on $A$.
We say that $v$ is a {\em Tate valuation} if the residual
valued field $(\kappa(v),\bar v)$ is a Tate valued field.
\end{enumerate}
\end{definition}

\begin{proposition}\label{prop_stays-valuation}
With the notation of definition {\em\ref{def_valuation-toplog}},
let $(K^{+\wedge},\cT_{K^+}^\wedge)$ (resp. $(K^\wedge,\cT^\wedge_K)$)
be the completion of the topological ring $(K^+,\cT_{K^+})$
(resp. $(K,\cT_K)$). We have :
\begin{enumerate}
\item
Let $\cT'$ be any topology on $K^+$ that is linear
and separated. Then either $\cT'=\cT_{K^+}$, or else
$\cT'$ is the discrete topology.
\item
$(K,|\cdot|_K)$ is a Tate valued field if and only if
$(K,\cT_K)$ is a Tate topological ring.
\item
$K^{+\wedge}$ is a valuation ring and $\cT_{K^+}^\wedge$
is the valuation topology on $K^{+\wedge}$.
\item
$K^\wedge$ is the field of fractions of $K^+$, and
$\cT_K^\wedge$ is the valuation topology on $K^\wedge$.
\item
The valuation $|\cdot|_K:K\to\Gamma_{\!\circ}$ extends
uniquely to a valuation $|\cdot|_K^\wedge:K^\wedge\to\Gamma_{\!\circ}$
whose valuation ring is $K^{+\wedge}$.
\end{enumerate}
\end{proposition}
\begin{proof}(i): Let $x\in K^+\setminus\{0\}$ be any element;
since $\cT'$ is separated, there exists an ideal $I\subset K^+$
that is open for the topology $\cT'$, and such that
$x\notin I$. Then $I\subset U_\gamma$, with $\gamma:=|x|_K$
(notation of definition \ref{def_valuation-toplog}(i)),
so $\cT'$ is finer than $\cT_{K^+}$. If $\cT'$ is not
the discrete topology, there exists $y\in I\setminus\{0\}$,
and therefore $U_\delta\subset I$, with $\delta:=|y|_K$,
so in this case $\cT'=\cT_{K^+}$.

(ii): Suppose that $\cT_K$ is f-adic; then, since $K^+$ is
a bounded subset of $K$, the valuation topology $\cT_{K^+}$
must be adic with a finitely generated ideal of adic
definition $I$ (proposition \ref{prop_f-adics}(ii)).
Such $I$ is necessarily principal, and $I\neq 0$, if
$\cT_K$ is not discrete; in this case, any generator
of $I$ is a topologically nilpotent unit in $K$.
Conversely, if $a\in K$ is a topologically nilpotent
unit for the valuation topology of $K$, then $a\in K^+$,
and it is easily seen that the valuation topology of
$K^+$ agrees with the $K^+a$-adic topology, so $\cT_{K^+}$
is adic and $\cT_K$ is f-adic.

(iii): Let $a_\bullet:=(a_i~|~i\in I)$ be any Cauchy
net in $K$ (indexed by a filtered ordered set $I$); for
every $\gamma\in\Gamma$, pick $i(\gamma)\in I$ such that
$|a_j-a_k|_K<\gamma$ for every $j,k\geq i(\gamma)$.
Suppose first that $|a_{i(\gamma)}|_K\leq\gamma$ for every
$\gamma\in\Gamma$; then $|a_j|_K\leq\gamma$ for every
$j\geq i(\gamma)$ and every $\gamma\in\Gamma$, {\em i.e.}
$a_\bullet$ converges to $0$, and we set
$|a_\bullet|^\wedge_K:=0$. Otherwise, there exists
$\delta\in\Gamma$ such that $|a_{i(\delta)}|_K>\delta$ and
then $|a_j|_K=\delta$ for every $j\geq i(\delta)$, in
which case we set $|a_\bullet|^\wedge_K:=|a_{i(\delta)}|_K$.
It is easily seen that $|a_\bullet|^\wedge_K$ depends only
on the class of $a_\bullet$ in $K^\wedge$, and if $b_\bullet$
is any other Cauchy sequence, then
$$
|(a_ib_i~|~i\in I)|^\wedge_K=
|a_\bullet|_K^\wedge\cdot|b_\bullet|^\wedge_K
\qquad\text{and}\qquad
|(a_i+b_i~|~i\in I)|^\wedge_K\leq
\max(|a_\bullet|_K^\wedge,|b_\bullet|^\wedge_K)
$$
(details left to the reader). It follows already that
$K^\wedge$ is a domain. Moreover, we have
$|a_\bullet|^\wedge_K\leq 1$ if and only if the class
of $a_\bullet$ is equivalent to the class of a Cauchy
net in $K^+$. Next, suppose that
$1\geq |a_\bullet|^\wedge_K\geq|b_\bullet|_K^\wedge>0$, in which
case there exists $i_0\in I$ such that
$$
|a_j|_K\geq|b_j|_K
\qquad\text{and}\qquad
|a_j|_K=\beta:=|a_\bullet|^\wedge_K
\qquad
\text{for every $j\geq i_0$}
$$
and set $v_j:=1$ for every $j\in I$ such that $j<i_0$
and $v_j:=a_j^{-1}b_j$ for every $j\geq i_0$. We claim that
$v_\bullet:=(v_j~|~j\in I)$ is a Cauchy net in $K^+$ :
indeed, for any $\gamma\in\Gamma$, pick $i(\gamma)\in I$
such that $|a_j-a_k|_K,|b_j-b_k|_R<\gamma$ for every
$j,k\geq i(\gamma)$; a simple calculation shows that
$|v_j-v_k|<\gamma\cdot\beta^2$ whenever
$j,k\geq\max(i_0,i(\gamma))$, as required. Now, since
$(v_ja_j~|~j\in I)$ and $b_\bullet$ represent the same
element in $K^{+\wedge}$, it follows that the class of
$a_\bullet$ divides the class of $b_\bullet$ in $K^{+\wedge}$
(the quotient is the class of $v_\bullet$). Lastly, (i)
easily implies that $\cT^\wedge_{K^+}$ is the valuation
topology on $K^{+\wedge}$, whence the assertion.

(v): Furthermore, notice that if the class $x$ of
$a_\bullet$ in $K^\wedge$ is not $0$, we have
$x=a_{i(\delta)}\cdot y$, with $y$ an invertible element
of $K^{+\wedge}$, so we must have $|x|'=|a_{i(\delta)}|_K$
for every valuation $|\cdot|'$ on $K^\wedge$ that extends
$|\cdot|_K$.

(iv): By the same token, we have
$x^{-1}=a_{i(\delta)}^{-1}\cdot y^{-1}$, which shows that
$K$ is a field, and concludes the proof.
\end{proof}

\begin{proposition}\label{prop_krasner}
Let $(K,|\cdot|_K)$ be a valued field such that the valuation
ring $K^+$ is henselian, $K^\mathrm{s}$ a separable closure of
$K$, and $(K^\wedge,|\cdot|^\wedge_K,\cT^\wedge_K)$ the completion
of $(K,|\cdot|_K,\cT_K)$. Then :
\begin{enumerate}
\item
The valuation ring $K^{\wedge+}$ of $K^\wedge$ is henselian.
\item
The ring $K^{\wedge\mathrm{s}}:=K^\wedge\otimes_KK^\mathrm{s}$
is an algebraic separably closed field extension of
$K^\wedge$.
\item
Especially, we have a natural group isomorphism :
$$
\Gal(K^\mathrm{s}/K)\isom\Gal(K^{\wedge\mathrm{s}}/K^\wedge)
\qquad
\sigma\mapsto K^\wedge\otimes_K\sigma.
$$
\end{enumerate}
\end{proposition}
\begin{proof}(i): The assertion follows immediately from
the more general :

\begin{claim} Let
$(A_\bullet,J_\bullet):=(A_\lambda,J_\lambda)~|~\lambda\in\Ob(\Lambda))$
be a system of henselian pairs indexed by any small category,
with transition morphisms $f_\phi:A_\lambda\to A_\mu$ such that
$f_\phi(J_\lambda)\subset J_\mu$ for every morphism
$\phi:\lambda\to\mu$ of $\Lambda$. Let $A$ (resp. $J$) be the
limit of the system $A_\bullet$ (resp. of the system $J_\bullet$);
then the pair $(A,J)$ is henselian.
\end{claim}
\begin{pfclaim} Let $\cA$ be the category whose objects are
all the pairs $(B,I)$, where $B$ is a ring and $I\subset B$ an
ideal, and whose morphisms $(B,I)\to(B',I')$ are all the ring
homomorphisms $f:B\to B'$ with $f(I)\subset I'$. Denote by
$\cA^\he$ the full subcategory of $\cA$ whose objects are the
henselian pairs. Then the inclusion functor $\cA^\he\to\cA$
admits a left adjoint, which assigns to every pair $(B,I)$
its henselization $(B^\he,I^\he)$. It is easily seen that $\cA$
is complete, so the claim follows directly from corollary
\ref{cor_mitchell}.
\end{pfclaim}

(ii): Let $K^\mathrm{a}$ be an algebraic closure of $K$ containing
$K^\mathrm{s}$ and contained in an algebraic closure
$K^{\wedge\mathrm{a}}$ of $K^\wedge$, and notice that the valuation
$|\cdot|_K$ extends uniquely to a valuation $|\cdot|_E$ on every
subextension $E\subset K^\mathrm{a}$, since $K^+$ is henselian
(\cite[Rem.6.1.12(iv)]{Ga-Ra}). Especially, for every such $E$
and every $K$-algebra homomorphism $\sigma:E\to K^\mathrm{a}$ we
must have
$$
|a|_E=|\sigma(a)|_{K^\mathrm{a}}
\qquad
\text{for every $a\in E$}.
$$

\begin{claim}\label{cl_Krasner}(Krasner's lemma)
With the foregoing notation, let $a\in K^\mathrm{s}$,
$b\in K^\mathrm{a}$, and suppose that for every homomorphism
$\sigma:K[a]\to K^\mathrm{s}$ of $K$-algebras with
$\sigma(a)\neq a$ we have
$$
|b-a|_{K^\mathrm{s}}<|\sigma(a)-a|_{K^\mathrm{s}}.
$$
Then $K[a]\subset K[b]$.
\end{claim}
\begin{pfclaim} It suffices to show that for every $K[b]$-algebra
homomorphism $\tau:K[a,b]\to K^\mathrm{a}$ we have $\tau(a)=a$.
However, suppose that $\tau(a)\neq a$ for such a map $\tau$;
we deduce that
$$
|b-a|_{K^\mathrm{a}}<|\tau(a)-a|_{K^\mathrm{s}}
\qquad\text{and}\qquad
|b-\tau(a)|_{K^\mathrm{a}}=|\tau(b-a)|_{K^\mathrm{a}}=
|b-a|_{K^\mathrm{a}}
$$
whence :
$$
|\tau(a)-a|_{K^\mathrm{s}}=|\tau(a)-b+b-a|_{K^\mathrm{s}}
\leq\max(|\tau(a)-b|_{K^\mathrm{a}},|b-a|_{K^\mathrm{a}})<
|\tau(a)-a|_{K^\mathrm{s}}
$$
a contradiction.
\end{pfclaim}

Let now $K^\wedge\subset E^\wedge$ be a non-trivial finite
Galois extension with $E^\wedge\subset K^{\wedge\mathrm{a}}$,
and pick $a\in E^\wedge$ such that $E^\wedge=K^\wedge[a]$. Let
$P(X)\in K^\wedge[X]$ be the minimal polynomial of $a$ over
$K^\wedge$; in light of (i) and \cite[Rem.6.1.12(iv)]{Ga-Ra},
the valuation $|\cdot|^\wedge_K$ of $K^\wedge$ extends uniquely
to a valuation of $|\cdot|^\wedge_E$ of $E^\wedge$; with this
notation, we set
$$
G:=\Gal(E^\wedge/K^\wedge)
\qquad
\gamma:=
\min(|\sigma(a)-a|^\wedge_E~|~\sigma\in G\setminus\one_{E^\wedge})
\in|E^\wedge|^\wedge_E.
$$
Say that $P=X^n+c_1X^{n-1}+\cdots+c_n$ and pick a polynomial
$Q(X)\in K[X]$ of degree $n$, say $Q=X^n+d_1X^{n-1}+\cdots+d_n$
such that
$$
|c_i-d_i|^\wedge_K<\gamma^n/|a|_E^{n-i}
\qquad
\text{for $i=1,\dots,n$}.
$$
It follows easily that $|Q(a)|^\wedge_E<\gamma^n$. Then,
say that $Q(X)=\prod_{i=1}^n(X-b_i)$ for some
$b_1,\dots,b_n\in K^\mathrm{a}$; we deduce that $|a-b|<\gamma$
for some $b\in\{b_1,\dots,b_n\}$, in which case claim
\ref{cl_Krasner} implies that $E^\wedge\subset K^\wedge[b]$,
and since $P$ and $Q$ have the same degree, we see that
$Q$ is irreducible over $K^\wedge$ and $E^\wedge=K^\wedge[b]$.
Especially, $b$ is separable over $K^\wedge$, so the roots
of $Q$ are all distinct, and therefore $b$ is also separable
over $K$. Summing up, we have already shown that
$K^\wedge K^\mathrm{s}\subset K^{\wedge\mathrm{a}}$ is a separable
closure of $K^\wedge$. It remains to check that the natural
map $K^\wedge\otimes_KK^\mathrm{s}\to K^\wedge K^\mathrm{s}$ is
injective. To this aim, let $K\subset E$ be any finite
Galois extension with $E\subset K^\mathrm{s}$, and denote by
$\tr_{E/K}:E\to K$ the $K$-linear trace map, so that
$\tr_{E/K}(a)=\sum_{\sigma\in\Gal(E/K)}\sigma(a)$.
Clearly $\tr_{E/K}$ is continuous for the valuation topologies
$\cT_E$ and $\cT_K$ of $E$ and $K$; on the other hand, since
$E$ is separable over $K$, the trace map induces a perfect
pairing $E\otimes_KE\to K$ : $x\otimes y\mapsto\tr_{E/K}(xy)$
for every $x,y\in E$. Hence, let $x_1,\dots,x_n$ be a basis of
the $K$-vector space $E$, and $y_1,\dots,y_n\in E$ the unique
elements such that $\tr_{E/K}(x_iy_j)=\delta_{ij}$ for every
$i,j=1,\dots,n$. There follow $K$-linear isomorphisms
$$
K^{\oplus n}\xrightarrow{\ \phi\ }E
\xrightarrow{\ \psi\ }K^{\oplus n}
\qquad\text{such that}\qquad
\psi\circ\phi=\one_{K^{\oplus n}}
$$
where $\phi(a_1,\dots,a_n):=\sum_{i=1}^na_ix_i$ and
$\psi(b):=(\tr_{E/K}(by_1),\dots,\tr_{E/K}(by_n))$ for every
$b\in E$ and every $(a_1,\dots,a_n)\in K^{\oplus n}$. Clearly
both $\phi$ and $\psi$ are continuous for the topology
$\cT_E$ and the product topology on $K^{\oplus n}$, hence
they are both homeomorphisms. Let $(E^\wedge,\cT^\wedge_E)$
be the completion of $(E,\cT_E)$; we deduce a commutative
diagram of $K^\wedge$-linear maps
$$
\xymatrix{
K^\wedge\otimes_KE \ar[rr]^-{K^\wedge\otimes_K\psi} \ar[d]_\mu & &
K^\wedge\otimes_KK^{\oplus n} \ar[d] \\
E^\wedge \ar[rr]^-{\psi^\wedge} & & (K^\wedge)^{\oplus n}
}$$
where $\psi^\wedge$ is the completion of $\psi$ and $\mu$
is the multiplication map; especially both horizontal arrows
are bijections, and the same holds for the right vertical
arrow. Then $\mu$ is bijective as well, whence the contention.

(iii) is an immediate consequence of (ii).
\end{proof}

\begin{remark} The left adjoint property of the topological
henselization $A\mapsto A^\he$ implies that for every f-adic
ring $A$, the identity map $\one_{A^\he}$ is the unique
continuous $A$-algebra endomorphism of $A^\he$. However,
$A^\he$ in general will admit also non-continuous $A$-algebra
endomorphisms. For instance, let $K$ be an algebraically
closed field of characteristic zero, and set
$$
F_0:=K[\Q_+]=K[T^{1/n}~|~n\in\N]
\qquad
F:=\Frac\,F_0.
$$
Let also $v:F_0\to\Q_\circ$ be the valuation of $F_0$
associated with the inclusion map $\Q_+\to\Q$ as in example
\ref{ex_toric-valuations} (where $\Q$ is endowed with its
standard ordering). Then $F^+:=\kappa(v)^+$ is a valuation
ring of $F$, and $(F,v)$ is a Tate valued field. Clearly
the value group of $v$ is $\Q$ and the residue field is
$K$. Let $\fm_F\subset F^+$ be the maximal ideal, and
$(F^{+\he},\fm_F^\he)$ the henselization of the pair
$(F^+,\fm)$; then $F^{+\he}$ is a valuation ring with
the same value group $\Q$, and the same residue field
$K$ (\cite[Lemma 6.2.5]{Ga-Ra}).
Moreover, the topological henselization $F^\he$ of the
Tate valued field $F$ is the field of fractions of
$F^{+\he}$. It then follows easily from
\cite[Prop.6.2.12]{Ga-Ra} that $F^\he$ is algebraically
closed; on the other hand, it is easily seen that $F$ is
not algebraically closed, so there exists a non-trivial
automorphism of $F^\he$ that fixes $F$.
\end{remark}

\begin{theorem}\label{th_val-cornerstone}
Let $K$ be a field, $A\subset K$ a subring, and
$\fp\subset A$ a prime ideal of $A$. Then there
exists a valuation ring $(V,\fm_V)$ of $K$ such that
$$
A\subset V
\qquad\text{and}\qquad
\fp=\fm_V\cap A.
$$
\end{theorem}
\begin{proof} We may replace $A$ by $A_\fp$, and assume
from start that $A$ is a local ring. Denote by $\cF$ the
set of all subrings $B$ of $K$ containing $A$ and such
that $1\notin\fp B$; we endow $\cF$ with the partial
ordering given by inclusion of subrings. Then $A\in\cF$,
and if $\cG$ is a totally ordered subset of $\cF$, then
$\bigcup_{B\in\cG}B$ still lies in $\cF$. By Zorn's lemma,
it follows that $\cF$ admits a maximal element $V$. First,
we claim that $V$ is local. Indeed, by assumption
$1\notin\fp V$, so $\fp V$ is contained in a maximal ideal
$\fm$ of $V$, and $1\notin\fp V_\fm$, so $V_\fm\in\cF$, and
therefore $V=V_\fm$, since $V$ is maximal. Let $\fm_V$ be
the maximal ideal of $V$; by construction,
$\fp\subset\fm_V\cap A$, and since $\fp$ is maximal in $A$,
we get $\fp=\fm_V\cap A$. It remains only to check that $V$
is a valuation ring of $K$. Thus, say that $x\in K$ and
$x\notin V$; by maximality of $V$, we have $V[x]\notin\cF$,
so $1\in\fp V[x]$, {\em i.e.} there exists an identity of
the form $1=a_0+a_1x+\cdots+a_nx^n$ for some
$a_0,\dots,a_n\in\fp V\subset\fm_V$. Since $1-a_0\in V^\times$,
we deduce a relation of the form
\set\begin{equation}\label{eq_for-x}
1=b_1x+\cdots+b_nx^n
\qquad
\text{where $b_1,\dots,b_n\in\fm_V$}.
\end{equation}
Likewise, if $x^{-1}\notin V$, we find an identity of the form
\set\begin{equation}\label{eq_for-x-inverse}
1=c_1x^{-1}+\cdots+c_mx^{-m}
\qquad
\text{where $c_1,\dots,c_m\in\fm_V$}
\end{equation}
and we may assume that $m,n\in\N$ are the minimal 
integers for which there exist identities \eqref{eq_for-x}
and \eqref{eq_for-x-inverse}. Up to swapping the roles
of $x$ and $x^{-1}$, we may also assume that $n\geq m$.
In this case we deduce
$$
1=b_1x+\cdots+b_{n-1}x^{n-1}+b_n\cdot(c_1x^{n-1}+\cdots+c_mx^{n-m})
$$
which is another identity of the type \eqref{eq_for-x},
but with a strictly smaller value of $n$, a contradiction.
We conclude that either $x\in V$ or $x^{-1}\in V$, as required. 
\end{proof}

\sset\subsubsection{}\label{subsec_dominant}
Let $K$ be any field, and denote by $\cL_K$ the set
of all local subrings of $K$. We endow $\cL_K$ with
a partial ordering, by declaring that, for every
$A,B\in\cL_K$ we have $A\leq B$ if and only if $B$
dominates $A$ (see definition \ref{def_dominant}(ii)).

\begin{corollary}\label{cor_cornerstone}
With the notation of \eqref{subsec_dominant}, the
following holds :
\begin{enumerate}
\item
For every $A\in\cL_K$ there exists a maximal element
of $\cL_K$ that dominates $A$.
\item
An element of $\cL_K$ is maximal if and only if it
is a valuation ring of $K$.
\end{enumerate}
\end{corollary}
\begin{proof}(i): Let $(B_i~|~i\in I)$ be a totally
ordered subset of $\cL_K$ such that $A\leq B_i$ for
every $i\in I$, and set $B:=\bigcup_{i\in I}B_i$.
It is easily seen that $B\in\cL_K$, and $A\leq B$.
Then the assertion follows from Zorn's lemma.

(ii): Let $A$ be any maximal element of $\cL_K$;
by theorem \ref{th_val-cornerstone} there exists
a valuation ring $V$ of $K$ that dominates $A$,
so $A=V$. Conversely, if $(V,\fm_V)$ is a valuation
ring of $K$ and $(A,\fm_A)$ any element of $\cL_K$
that dominates $V$, then we must have $A=V$ : indeed,
if the latter fails, there exists $a\in A\setminus V$,
so that $a^{-1}\in\fm_V\subset\fm_A$, which is absurd.
\end{proof}

\sset\subsubsection{}
Henceforth, and until the end of this section, we let
$(K,|\cdot|)$ be a given valued field, whose valuation
ring (resp. maximal ideal, resp. residue field, resp.
value group) shall be denoted $K^+$ (resp. $\fm_K$, resp.
$\kappa$, resp. $\Gamma$). The following result is due
to Nagata (\cite[Th.3]{Na2}).

\begin{proposition}\label{prop_Gruson-out}
Let $A$ be an essentially finitely presented $K^+$-algebra,
$M$ a finitely generated $K^+$-flat $A$-module. Then $M$ is
a finitely presented $A$-module.
\end{proposition}
\begin{proof} Let us write $A=S^{-1}B$ for some finitely
presented $K^+$-algebra $B$, and some multiplicative subset
$S\subset B$; then we may find a finitely generated $B$-module
$M_B$ with an isomorphism $\phi:S^{-1}M_B\isom M$ of $A$-modules.
Let $M'_B$ be the image of $M_B$ in $M_B\otimes_{K^+}K$;
since $M$ is $K^+$-flat, $\phi$ induces an isomorphism
$S^{-1}M'_B\to M$. It then suffices to show that $M'_B$
is a finitely presented $B$-module; hence we may replace
$A$ by $B$ and $M$ by $M'_B$, and assume that $A$ is a
finitely presented $K^+$-algebra.

We are further easily reduced to the case where
$A=K^+[T_1,\dots,T_n]$ is a free polynomial $K^+$-algebra.
Pick a finite system of generators $\Sigma$ for $M$; define
increasing filtrations $\Fil_\bullet A$ and $\Fil_\bullet M$
on $A$ and $M$, by letting $\Fil_kA$ be the $K^+$-submodule of
all polynomials $P(T_1,\dots,T_n)\in A$ of total degree
$\leq k$, and setting $\Fil_kM:=\Fil_kA\cdot\Sigma\subset M$
for every $k\in\N$. We consider the Rees algebra
$\sR(\underline A)_\bullet$ and the Rees module
$\sR(\underline M)_\bullet$ associated with these filtrations
as in definition \ref{def_Rees}(iii,iv). Say that $\Sigma$
is a subset of cardinality $N$; we obtain an $A$-linear
surjection $\phi:A^{\oplus N}\to M$, and we set
$M'_k:=(\Fil_kA)^{\oplus N}\cap\Ker\,\phi$ for every $k\in\N$.
Notice that the resulting graded $K^+$-module
$M'_\bullet:=\bigoplus_{k\in\N}M'_k$ is actually a
$\sR(\underline A)_\bullet$-module and we get a short exact
sequence of graded $\sR(\underline A)_\bullet$-modules
$$
C_\bullet
\quad :\quad
0\to M'_\bullet\to\sR(\underline A)_\bullet^{\oplus N}\to
\sR(\underline M)_\bullet\to 0.
$$
Notice also that the $K^+$-modules $\sR(\underline A)_k^{\oplus N}$
and $\sR(\underline M)_k$ are torsion-free and finitely
generated, hence they are free of finite rank (see
\cite[Rem.6.1.12(ii)]{Ga-Ra}); then the same holds for $M'_k$,
for every $k\in\N$. It follows that the complex
$C_\bullet\otimes_{K^+}\kappa$ is still short exact.
On the other hand, $\sR(\underline A)_\bullet$ is a $K^+$-algebra
of finite type (see example \ref{ex_Rees-free}), hence
$R:=\sR(\underline A)_\bullet\otimes_{K^+}\kappa$ is a noetherian
ring; therefore, $M'_\bullet\otimes_{K^+}\kappa$ is a graded
$R$-module of finite type, say generated by the system of
homogenenous elements $\{\bar x_1,\dots,\bar x_t\}$. Lift these
elements to a system $\Sigma':=\{x_1,\dots,x_t\}$ of homogeneous
elements of $M'_\bullet$, and denote by
$M''_\bullet\subset M'_\bullet$ the
$\sR(\underline A)_\bullet$-submodule generated by $\Sigma'$.
By construction, $M''_k\otimes_{K^+}\kappa=M'_k\otimes_{K^+}\kappa$,
hence $M''_k=M'_k$ for every $k\in\N$, by Nakayama's lemma.
In other words, $\Sigma'$ is a system of generators for
$M'_\bullet$; it follows easily that the image of $\Sigma'$
in $\Ker\,\phi$ is also a system of generator for the latter
$A$-module, and the proposition follows.
\end{proof}

\begin{corollary}\label{cor_coherence}
Every essentially finitely presented $K^+$-algebra is a
coherent ring.
\end{corollary}
\begin{proof} One reduces easily to the case of a free
polynomial $K^+$-algebra $K^+[T_1,\dots,T_n]$, in which
case the assertion follows immediately from proposition
\ref{prop_Gruson-out} : the details shall be left to
the reader.
\end{proof}

\begin{lemma}\label{lem_fin-min}
If $A$ if a finitely presented $K^+$-algebra, the subset
$\mathrm{Min}\,A\subset\Spec\,A$ of minimal prime ideals
is finite.
\end{lemma}
\begin{proof} We begin with the following :
\begin{claim}\label{cl_going-down}
The assertion holds when $A$ is a flat $K^+$-algebra.
\end{claim}
\begin{pfclaim} Indeed, in this case, by the going-down
theorem (\cite[Th.9.5]{Mat}) the minimal primes lie in
$\Spec\,A\otimes_{K^+}K$, which is a noetherian ring.
\end{pfclaim}

In general, $A$ is defined over a $\Z$-subalgebra $R\subset K^+$
of finite type. Let $L\subset K$ be the field of fractions of
$R$, and set $L^+:=L\cap K^+$; then $L^+$ is a valuation ring
of finite rank, and there exists a finitely presented $L^+$-algebra
$A'$ such that $A\simeq A'\otimes_{L^+}K^+$.
For every prime ideal $\fp\subset L^+$, let $\kappa(\fp)$
be the residue field of $L^+_\fp$, and set :
$$
A'(\fp):=A'\otimes_{L^+}\kappa(\fp)
\qquad
K^+(\fp):=K^+\otimes_{L^+}\kappa(\fp).
$$
Since $\Spec\,L^+$ is a finite set, it suffices to show that
the subset $\mathrm{Min}\,A\otimes_{L^+}\kappa(\fp)$ is finite
for every $\fp\in\Spec\,L^+$. However,
$A\otimes_{L^+}\kappa(\fp)\simeq A'(\fp)\otimes_{\kappa(\fp)}K^+(\fp)$,
so this ring is a flat $K^+(\fp)$-algebra of finite presentation.
Let $I$ be the nilradical of $K^+(\fp)$; since $K^+(\fp)/I$ is a
valuation ring, claim \ref{cl_going-down} applies with $A$ replaced
by $A\otimes_{K^+}K^+(\fp)/I$, and concludes the proof.
\end{proof}

\begin{lemma}\label{lem_integral-fp-ext}
Let $(A,\fm_A)$ be a local ring, $K^+\to A$ a local and
essentially finitely presented ring homomorphism.
Then there exists a local morphism $\phi:V\to A$ of essentially
finitely presented $K^+$-algebras, such that the following holds :
\begin{enumerate}
\item
$V$ is a valuation ring, its maximal ideal $\fm_V$ is generated
by the image of $\fm_K$, and the map $K^+\to V$ induces
an isomorphism $\Gamma\isom\Gamma_{\!V}$ of value groups.
\item
$\phi$ induces a finite field extension $V/\fm_V\to A/\fm_A$.
\end{enumerate} 
\end{lemma}
\begin{proof} Set $X:=\Spec\,A$, $S:=\Spec\,K^+$ and let
$x\in X$ be the closed point. Pick $a_1,\dots,a_d\in A$ whose
classes in $\kappa(x)$ form a transcendence basis over $\kappa(s)$.
The system $(a_i~|~i\leq d)$ defines a factorization of the
morphism $\Spec\,\phi:X\to S$ as a composition
$X\xrightarrow{g}Y:=\A^d_S\xrightarrow{h}S$, such that
$\xi:=g(x)$ is the generic point of $h^{-1}(s)\subset Y$.
The morphism $g$ is essentially finitely presented
(\cite[Ch.IV, Prop.1.4.3(v)]{EGAIV}) and moreover :

\begin{claim}\label{cl_it-s-a-Gauss}
The stalk $\cO_{Y,\xi}$ is a valuation ring with value
group $\Gamma$.
\end{claim}
\begin{pfclaim} Indeed, set $B:=K^+[T_1,\dots,T_d]$; one sees
easily that $V:=\cO_{\!X,\xi}$ is the valuation ring of the
Gauss valuation $|\cdot|_B:\Frac(B)\to\Gamma\cup\{0\}$ such that
$$
\Bigl|\sum_{\alpha\in\N^d}a_\alpha T^\alpha\Bigr|_B=
\max\{|a_\alpha|~|~\alpha\in\N^d\}
$$
(where $a_\alpha\in K^+$ and
$T^\alpha:=T_1^{\alpha_1}\cdots T_r^{\alpha_r}$ for all
$\alpha:=(\alpha_1,\dots,\alpha_d)\in\N^d$, and $a_\alpha=0$
except for finitely many $\alpha\in\N^d$).
\end{pfclaim}

In view of claim \ref{cl_it-s-a-Gauss}, to conclude the proof 
it suffices to notice that $\kappa(x)$ is a finite extension of
$\kappa(\xi)$.
\end{proof}

\begin{proposition}\label{prop_integral-fp-ext}
Let $\phi:K^+\to A$ and $\psi:A\to B$ be two essentially
finitely presented ring homomorphisms. Then :
\begin{enumerate}
\item
If\/ $\psi$ is integral, $B$ is a finitely presented $A$-module.
\item
If $A$ is local, $\phi$ is local and flat, and $A/\fm_K A$ is
a field, then $A$ is a valuation ring and $\phi$ induces an
isomorphism $\Gamma\isom\Gamma_{\!\!A}$ of value groups.
\end{enumerate}
\end{proposition}
\begin{proof} (i): The assumption means that there exist a finitely
presented $A$-algebra $C$ and a multiplicative system $T\subset C$
such that $B=T^{-1}C$. Let $I:=\bigcup_{t\in T}\Ann_C(t)$. Then
$I$ is the kernel of the localization map $C\to B$;
the latter is integral by hypothesis, hence for every $t\in T$
there is a monic polynomial $P(X)\in C[X]$, say of degree $n$,
such that $P(t^{-1})=0$ in $B$, hence $t't^n\cdot P(t^{-1})=0$
in $C$, for some $t'\in T$, {\em i.e.} $t'(1-ct)=0$ holds in $C$
for some $c\in C$, in other words, the image of $t$ is
already a unit in $C/I$, so that $B=C/I$. Then lemma
\ref{lem_from-I-to-Z} says that $V(I)\subset\Spec\,C$ is closed under
generizations, hence $V(I)$ is open, by corollary \ref{cor_min-max}
and lemma \ref{lem_fin-min}, and finally $I$ is finitely generated,
by lemma \ref{lem_flat-quot}(ii). So $B$ is a finitely presented
$A$-algebra. Then (i) follows from the well known :

\begin{claim}\label{cl_fin-present}
Let $R$ be any ring, $S$ a finitely presented and integral
$R$-algebra. Then $S$ is a finitely presented $R$-module.
\end{claim}
\begin{pfclaim} Let us pick a presentation :
$S=R[x_1,\dots,x_n]/(f_1,\dots,f_m)$. By assumption,
for every $i\leq n$ we can find a monic polynomial
$P_i(T)\in R[T]$ such that $P_i(x_i)=0$ in $S$.
Let us set $S':=R[x_1,\dots,x_n]/(P_1(x_1),\dots,P_n(x_n))$;
then $S'$ is a free $R$-module of finite rank, and there
is an obvious surjection $S'\to S$ of $R$-algebras, the
kernel of which is generated by the images of the
polynomials $f_1,\dots, f_m$. The claim follows.
\end{pfclaim}

(ii): By lemma \ref{lem_integral-fp-ext}, we may assume from
start that the residue field of $A$ is a finite extension of
$\kappa$. We can write $A=C_\fp$ for a $K^+$-algebra $C$ of
finite presentation, and a prime ideal $\fp\subset C$
containing $\fm_K$; under the stated assumptions, $A/\fm_K A$
is a finite field extension of $\kappa$, hence $\fp$ is a
maximal ideal of $C$, and moreover $\fp$ is isolated in the
fibre $g^{-1}(s)$ of the structure morphism $g:Z:=\Spec\,C\to\Spec\,K^+$.
It then follows from \cite[Ch.IV, Cor.13.1.4]{EGAIV-3} that,
up to shrinking $Z$, we may assume that $g$ is a quasi-finite
morphism (\cite[Ch.II, \S6.2]{EGAII}). Let $K^{\he+}$ be the
henselization of $K^+$, and set $C^\he:=C\otimes_{K^+}K^{\he+}$.
There exists precisely one prime ideal $\fq\subset C^\he$ lying
over $\fp$, and $C^\he_\fq$ is a flat, finite and finitely
presented $K^{\he+}$-algebra by \cite[Ch.IV, Th.18.5.11]{EGA4},
so $C^\he_\fq$ is henselian, by \cite[Ch.IV, Prop.18.5.6(i)]{EGA4},
and then $C^\he_\fq$ is the henselization of $A$, since the
natural map $A\to C^\he_\fq$ is ind-{\'e}tale and induces an
isomorphism $C^\he_\fq\otimes_{K^+}\kappa\simeq\kappa(x)$.
Finally, $C^\he_\fq$ is also a valuation ring with value
group $\Gamma$, by virtue of
\cite[Lemma 6.1.13 and Rem.6.1.12(vi)]{Ga-Ra}. To conclude
the proof of (ii), it now suffices to remark :

\begin{claim}\label{cl_finally-named}
Let $V$ be a valuation ring, $\phi:R\to V$ a faithfully
flat morphism. Then $R$ is a valuation ring and $\phi$
induces an injection
$\phi_\Gamma:\Gamma_{\!\!R}\to\Gamma_{\!\!V}$ of the
respective value groups.
\end{claim}
\begin{pfclaim}[] Clearly $R$ is a domain, since it is
a subring of $V$. To show that $R$ is a valuation ring, it
then suffices to prove that, for any two ideals
$I,J\subset R$, we have either $I\subset J$ or $J\subset I$.
However, since $V$ is a valuation ring, we know already
that either $I\cdot V\subset J\cdot V$ or
$J\cdot V\subset I\cdot V$; since $\phi$ is faithfully flat,
the assertion follows. By the same token, we deduce that
an ideal $I\subset R$ is equal to $R$ if and only if
$I\cdot V=V$, which implies that $\phi_\Gamma$ is injective
(see \cite[\S6.1.11]{Ga-Ra}).
\end{pfclaim}
\end{proof}

\subsection{Huber's theory of the valuation spectrum}
\label{sec_Spv-of-ring}
In this section we present the first elements of Huber's
theory of the valuation spectrum of an arbitrary ring, for
which the original reference is his Habilitationschrift
\cite{Hu1}. Let $A$ be any ring. We denote by
$$
\Spv\,A
$$
the set of equivalence classes of valuations on $A$.
For every $a,b\in A$ we set
$$
R_A\Bigl(\frac{a}{b}\Bigr):=\{v\in\Spv\,A~|~v(a)\leq v(b)\neq 0\}
$$
and more generally, with every $n\in\N$ and every
$a_1,\dots,a_n,b_1,\dots,b_n\in A$ we associate the
{\em rational subset}
$$
R_A\Bigl(\frac{a_1}{b_1},\dots,\frac{a_n}{b_n}\Bigr):=
\bigcap_{i=1}^nR_A\Bigl(\frac{a_i}{b_i}\Bigr)
$$
and we denote $\cB_A$ the boolean algebra generated
by the system of rational subsets of $\Spv\,A$.
Explicitly, every element of $\cB_A$ is a finite
union of subsets of the form
$$
\{v\in\Spv\,A~|~v(a_1)<v(b_1),\dots,v(a_n)<v(b_n),\ 
                v(c_1)\leq v(d_1),\dots,v(c_m)\leq v(d_m)\}
$$
where $n,m\in\N$ are arbitrary integers, and
$a_1,b_1,\dots,a_n,b_n,c_1,d_1,\dots,c_m,d_m\in A$ are
arbitrary elements.

\begin{definition} For every ring $A$, we endow $\Spv\,A$
with the topology $\cT_A$ with basis given by the rational
subsets. The pair $(\Spv\,A,\cT_A)$ will be called the
{\em valuation spectrum} of $A$.
\end{definition}

\begin{theorem}\label{th_Spv-spectral}
For every ring $A$, the topological space $(\Spv\,A,\cT_A)$
is spectral, and $\cB_A$ is the set of all constructible
subsets of\/ $\Spv\,A$.
\end{theorem}
\begin{proof} We endow $\Spv\,A$ with the topology
$\cT_\mathrm{cons}$ with basis $\cB_A$. Endow as well the
set $\{0,1\}$ with its discrete topology, and the set
$\cP:=\{0,1\}^{A\times A}$ with the product topology.
We attach to every valuation $v$ of $A$ a mapping
$\phi_v:A\times A\to\{0,1\}$, by declaring that
$\phi_v(a,b)=1$ if and only if $v(a)\geq v(b)$.
If $v,w$ are any two valuations on $A$, then it is
clear that $\phi_v=\phi_w$ if and only if $v$ is
equivalent to $w$, so that the rule $v\mapsto\phi_v$
yields a well defined injective map
$$
\phi:\Spv\,A\to\cP
$$
and it is easily seen that the topology on $\Spv\,A$
induced by $\cP$ via $\phi$ agrees with $\cT_\mathrm{cons}$.

\begin{claim}\label{cl-correct-conditions}
The image of $\phi$ is the subset of $\cP$ of all mappings
$f:A\times A\to\{0,1\}$ fulfilling the following conditions
for every $a,b,c\in A$ :
\begin{enumerate}
\alphaenu
\item
we have either $f(a,b)=1$ or $f(b,a)=1$
\item
if $f(a,b)=f(b,c)=1$, then $f(a,c)=1$
\item
we have either $f(a,a+b)=1$ or $f(b,a+b)=1$
\item
if $f(a,b)=1$, then $f(ac,bc)=1$
\item
if $f(ac,bc)=1$ and $f(0,c)=0$, then $f(a,b)=1$
\item
$f(0,1)=0$.
\end{enumerate}
\end{claim}
\begin{pfclaim} It is clear that $\phi_v$ fulfills these
conditions for every $v\in\Spv\,A$. Conversely, let $f$
be a mapping fulfilling conditions (a)--(f); from (f),
(a) and (d) we deduce :
\begin{enumerate}
\addenu\addenu\addenu\addenu\addenu\addenu
\alphaenu
\item
$f(x,0)=1$ for every $x\in A$.
\end{enumerate}
Set $\fp:=\{a\in A~|~f(0,a)=1\}$. We show that $\fp$ is an
ideal of $A$. Indeed, from (a) it follows that $0\in\fp$.
Next, say that $a,b\in\fp$; from (c) we may assume that
$f(a,a+b)=1$, and then (b) implies that $a+b\in\fp$.
Lastly, if $a\in\fp$ and $b\in A$, then (d) implies that
$ab\in\fp$, as required. Moreover, (e) and (f) imply that
$\fp$ is a prime ideal, and we set $\kappa:=\Frac\,A/\fp$.
We notice :
\begin{enumerate}
\addenu\addenu\addenu\addenu\addenu\addenu\addenu
\alphaenu
\item
$f(a+x,a)=f(a,a+x)=1$ for every $a\in A$ and $x\in\fp$.
\end{enumerate}
Indeed, suppose that $f(a+x,a)=0$ for such $a$ and $x$; then
(c) implies that $f(-x,a)=1$, and combining with (b) we get
$f(0,a)=1$. On the other hand, (g) says that $f(a+x,0)=1$,
and combining again with (b) we derive $f(a+x,a)=1$, a
contradiction. To derive the second identity, set $b:=a+x$,
so that $f(b-x,b)=1$, by the foregoing, as required.

We can now show that $f$ factors through a mapping
$$
\bar f:A/\fp\times A/\fp\to\{0,1\}.
$$
Indeed, let $a,b\in A$ and $c\in\fp$ be arbitrary elements,
and assume first that $f(a,b)=1$; we need to show that
$f(a+c,b)=1$, and this follows from (b) and (h). Similarly,
we get $f(a,b+c)=1$. Lastly, if $f(a,b)=0$ we must have
$f(a+c,b)=0$, for otherwise the foregoing would give
$f(a,b)=f(a+c-c,b)=f(a+c,b)=1$, a contradiction; similarly,
we see that $f(a,b+c)=0$ in this case. This shows that
$\bar f$ is well defined, and by construction we have
\begin{enumerate}
\addenu\addenu\addenu\addenu\addenu\addenu\addenu\addenu
\alphaenu
\item
$\bar f(0,x)=0$ for every $x\in A/\fp\setminus\{0\}$.
\end{enumerate}
Now, for every $x\in\kappa$, write $x=a^{-1}b$ for some
$a,b\in A/\fp$, with $a\neq 0$; we say that $x$ is
{\em $f$-positive} if $\bar f(a,b)=1$. From (i) and (e)
it follows easily that this condition does not depend
on the choice of $a$ and $b$, and we let $V\subset\kappa$
be the subset of all $f$-positive elements. Suppose that
$x,y\in V$, and write $x=a^{-1}b$, $y=a^{-1}b'$ for some
$a,b,b'\in A/\fp$ with $a\neq 0$; this means that
$\bar f(a,b)=\bar f(a,b')=1$, and by virtue of (c) we
may assume that $\bar f(b,b+b')=1$. Then (b) implies that
$\bar f(a,b+b')=1$, so $x+y\in V$. Likewise, (d) implies
$\bar f(aa',ba')=1$ and $\bar f(ba',bb')=1$, and then
(b) yields $f(aa',bb')=1$, so $xy\in V$. Also, we have
either $\bar f(-1,1)=1$ or $\bar f(1,-1)=1$ by (a), so
$-1\in V$; summing up, we conclude that $V$ is a subring
of $\kappa$, and (a) implies that $V$ is a valuation
ring of $\kappa$. Now, let $w$ be the valuation of $\kappa$
attached to $V$ as in remark \ref{rem_valuations}(iv),
and set $v:=w\circ\pi$, where $\pi:A\to\kappa$ is the
natural map; it is easily seen that $\phi_v=f$, whence
the claim.
\end{pfclaim}

Claim \ref{cl-correct-conditions} implies that the
image of $\phi$ is a closed subset of the quasi-compact
and separated topological space $\cP$, so
$(\Spv\,A,\cT_\mathrm{cons})$ is quasi-compact and separated
as well. Now, it is clear that every rational subset of
$\Spv\,A$ is open and closed in the topology $\cT_\mathrm{cons}$.
Moreover, it follows easily from remark \ref{rem_valuations}(vi)
that $(\Spv\,A,\cT_A)$ is a $T_0$ topological space. Then the
theorem follows from proposition \ref{prop_T_0-criterion}.
\end{proof}

\begin{remark}\label{rem_Spv-of-ring}
(i)\ \
Let $A$ be any ring. Quite generally, we have
$$
R_A\Bigl(\frac{f_1}{f_0},\dots,\frac{f_n}{f_0}\Bigr)\cap
R_A\Bigl(\frac{g_1}{g_0},\dots,\frac{g_m}{g_0}\Bigr)=
R_A\Bigl(\frac{f_ig_j}{f_0g_0}~|~i=0,\dots,n,\ j=0,\dots,m\Bigr)
$$
for every sequences $f_\bullet:=(f_0,\dots,f_n)$
and $g_\bullet:=(g_0,\dots,g_m)$ of elements of $A$.

(ii)\ \
Every ring homomorphism $f:A\to B$ induces a mapping
$$
\Spv\,f:\Spv\,B\to\Spv\,A
\qquad
v\mapsto v\circ f
$$
and obviously we have
$$
(\Spv\,f)^{-1}R_A\Bigl(\frac{a}{b}\Bigr)=
R_B\Bigl(\frac{f(a)}{f(b)}\Bigr)
\qquad
\text{for every $a,b\in A$}
$$
so that $\Spv\,f$ is a continuous spectral map
$(\Spv\,B,\cT_B)\to(\Spv\,A,\cT_A)$ (remark
\ref{rem_sorite-qcoh-maps}(iii)).

(iii)\ \
Moreover, we have a natural {\em support map}
$$
\sigma_A:\Spv\,A\to\Spec\,A
\qquad
v\mapsto\Ker\,v
$$
and notice that
$$
\sigma_A^{-1}(\Spec\,A_s)=R_A\Bigl(\frac{s}{s}\Bigr)
\qquad
\text{for every $s\in A$}
$$
so also $\sigma_A$ is spectral. Furthermore, it is easily
seen that $\sigma_A$ restricts to a homeomorphism
$$
(\Spv\,A)_0\isom\Spec\,A
$$
where $(\Spv\,A)_0$ is the subset of all rank zero
valuations on $A$, endowed with the topology induced
by the inclusion map into $\Spv\,A$. Notice as well
that, by remark \ref{rem_semi-norm}(iv,v), we have a
natural homeomorphism
\set\begin{equation}\label{eq_fibers-of-supp}
\sigma_A^{-1}(\fp)\isom\Spv\,\kappa(\fp)
\qquad
\text{for every $\fp\in\Spec\,A$}
\end{equation}
where the fiber $\sigma^{-1}_A(\fp)$ is endowed with the
topology induced by $\Spv\,A$ via the inclusion map.

(iv)\ \
Let $S\subset A$ be any multiplicative system, and
$i:A\to S^{-1}A$ the localization map; taking into
account remark \ref{rem_semi-norm}(iv), it is easily
seen that $\Spv\,i$ is injective, and
$$
\Img\,\Spv\,i=\bigcap_{s\in S}R_A\Bigl(\frac{s}{s}\Bigr).
$$
Moreover, we have
$$
R_{S^{-1}A}\Bigl(\frac{s^{-1}a}{t^{-1}b}\Bigr)=
(\Spv\,i)^{-1}R_A\Bigl(\frac{at}{bs}\Bigr)
\qquad
\text{for every $a,b\in R$ and $s,t\in S$}
$$
so the topology of $\Spv\,S^{-1}A$ is induced by
the topology of $\Spv\,A$, via the mapping $\Spv\,i$,
and therefore the latter identifies $\Spv\,S^{-1}A$
with a pro-constructible subset of $\Spv\,A$.

(v)\ \
For any ring $A$, we set
$$
\Spv^+A:=
\{v\in\Spv\,A~|~\text{$v(a)\leq 1$ for every $a\in A$}\}.
$$
Hence, $\Spv^+A$ is a pro-constructible subset of
$\Spv\,A$, and we endow it with the topology induced
by the inclusion map into $\Spv\,A$, so that $\Spv^+A$
is a spectral space (theorem \ref{th_Spv-spectral} and
corollary \ref{cor_procon-is-spec}). For every $v\in\Spv^+A$,
let also
$$
\sigma_A^+(v):=\{a\in A~|~v(a)<1\}.
$$
It is easily seen that $\sigma_A^+(v)$ is a prime ideal
of $A$, so we get a well defined map
$$
\sigma_A^+:\Spv^+A\to\Spec\,A.
$$
Notice that $\sigma^+_A$ is surjective : indeed, for
every prime ideal $\fp$, the trivial valuation $v_\fp$
with support equal to $\fp$ lies in $\Spv^+A$, and
$\sigma^+_A(v_\fp)=\fp$. Moreover, we have
$$
(\sigma^+_A)^{-1}(\Spec\,A_s)=
R_A\Bigl(\frac{1}{s}\Bigr)\cap\Spv^+A
\qquad
\text{for every $s\in A$}
$$
so $\sigma^+_A$ is continuous and spectral (corollary
\ref{cor_procon-is-spec}). In many references, the prime
ideal $\sigma_A^+(v)$ is called the {\em center} of the
valuation $v$. Notice also that, for every ring homomorphism
$f:A\to B$, the map $\Spv\,f$ restricts to a continuous
spectral map
$$
\Spv^+f:\Spv^+B\to\Spv^+A.
$$
\end{remark}

\sset\subsubsection{}\label{subsec_Spv-of-scheme}
Let $X$ be any scheme. A {\em valuation of $X$} is a
pair $(x,v)$ consisting of a point $x\in X$ and the
equivalence class of a valuation $v$ of the residue
field $\kappa(x)$; then we say that $x$ is the
{\em support} of $(x,v)$. We denote by
$$
\Spv\,X
$$
the set of valuations of $X$, and we endow it with
the coarsest topology $\cT_X$ containing the subsets
$$
U\Bigl(\frac{a}{b}\Bigr):=
\{(x,v)\in\Spv\,X~|~x\in U,\ v(a(x))\leq v(b(x))\neq 0\}
$$
where $U$ ranges over all the open subsets of $X$, and
$a,b$ are arbitrary elements of $\cO_{\!X}(U)$ (and
$a(x),b(x)$ are the images of $a$ and $b$ in $\kappa(x)$).
The topological space $(\Spv\,X,\cT_X)$ is called the
{\em valuation spectrum of\/ $X$}. Any morphism of schemes
$f:Y\to X$ induces a mapping
$$
\Spv\,f:\Spv\,Y\to\Spv\,X
\qquad
(y,w)\mapsto(f(y),w\circ f^\flat_y)
$$
where $f^\flat_y:\kappa(f(y))\to\kappa(y)$ is the ring
homomorphism induced by the morphism of structure sheaves
$f^\flat:\cO_{\!X}\to f_*\cO_Y$ associated with $f$. Notice
that
$$
(\Spv f)^{-1}U\Bigl(\frac{a}{b}\Bigr)=
f^{-1}U\Bigl(\frac{f^\flat_U(a)}{f^\flat_U(b)}\Bigr)
$$
for every $U$, $a$ and $b$ as in the foregoing, so
$\Spv\,f$ is a continuous mapping. In the same vein,
consider the {\em support map}
$$
\sigma_X:\Spv\,X\to X
\qquad
(x,v)\to x.
$$
Clearly $\sigma_X^{-1}U=U(\frac{1}{1})$ for every open
subset $U\subset X$, so $\sigma_X$ is continuous as well.

\begin{remark}
The construction of \eqref{subsec_Spv-of-scheme}
generalizes that of \eqref{sec_Spv-of-ring} : indeed,
let $A$ be any ring; we have a natural mapping
\set\begin{equation}\label{eq_Spv-of-Spec}
\Spv\,A\to\Spv(\Spec\,A)
\qquad
v\mapsto(\Ker\,v,\bar v)
\end{equation}
where, for any valuation $v$ on $A$, we let $\bar v$ be
the residual valuation of $v$ (see remark
\ref{rem_semi-norm}(v)). In light of remarks
\ref{rem_valuations}(vi) and \ref{rem_Spv-of-ring}(iv)
it is easily seen that \eqref{eq_Spv-of-Spec} is
a homeomorphism (details left to the reader). Taking into
account theorem \ref{th_Spv-spectral}, we see that $\Spv\,X$
is locally spectral for every scheme $X$, and if $X$ is
quasi-compact and quasi-separated, then $\Spv\,X$ is
spectral (lemma \ref{lem_charact-coh-sob-sp}(iv)).
Moreover, remark \ref{rem_Spv-of-ring}(ii) and proposition
\ref{prop_local-qc-qs}(iii) imply that $\Spv\,f$ is a
spectral map, for every morphism of schemes $f:Y\to X$.
\end{remark}

\begin{proposition}\label{prop_tensor-Spv}
Let $f:A\to B$ and $g:A\to A'$ be two ring homomorphisms,
and set $B':=A'\otimes_AB$. Then the induced continuous
map of topological spaces
$$
\phi:\Spv\,B'\to\Spv\,A'\times_{\Spv A}\Spv\,B
$$
is surjective.
\end{proposition}
\begin{proof} Let $v'\in\Spv\,A'$ and $w\in\Spv\,B$ be two
elements such that $v:=\Spv(g)(v')=\Spv(f)(w)$, and denote by
$$
\bar v\in\Spv\,\kappa(v)
\qquad
\bar v{}'\in\Spv\,\kappa(v')
\qquad
\bar w\in\Spv\,\kappa(w)
$$
the residual valuations of $v$, $v'$, and $w$; the maps $f$
and $g$ induce field extensions $\kappa(v)\to\kappa(v')$
and $\kappa(v)\to\kappa(w)$, whence a ring homomorphism
$h:B'\to C:=\kappa(v')\otimes_{\kappa(v)}\kappa(w)$.
Suppose that $u\in\Spv\,C$ is an element whose image in
$\Spv\,\kappa(v')\times_{\Spv\,\kappa(v)}\Spv\,\kappa(w)$
equals the pair $(\bar v{}',\bar w)$; then
$\phi\circ\Spv(h)(u)=(v',w)$. Thus, we are reduced to
the case where $A$, $A'$ and $B$ are three fields.
To ease notation, set $A^+:=\kappa(v)^+$, and define
likewise $A'^+$ and $B^+$; let as well
$B'^+:=A'^+\otimes_{A^+}B^+$. Since the maps $A^+\to A'^+$
and $A^+\to B^+$ are flat, it is easily seen that the
induced map $B'^+\to B'$ is injective. Pick any prime
ideal $\fp\in\Spec\,B'^+$ whose image in $\Spec\,A'^+$
(resp. in $\Spec\,B^+$) is the (unique) closed point;
the induced map $B'^+_\fp\to B'_\fp$ is still injective,
and especially, $B'_\fp\neq 0$, so we may find a
maximal ideal $\fm$ of $B'_\fp$, and we denote by
$\kappa(\fm)$ the residue field of $\fm$. The image
of $B'^+_\fp$ in $\kappa(\fm)$ is a local ring $D$, and
by corollary \ref{cor_cornerstone} we may find a
valuation ring $(V,\fm_V)$ of $\kappa(\fm)$ that
dominates $D$. Lastly, $V$ yields a well defined point
$u\in\Spv\,B'$, and a simple inspection shows that
$\phi(u)=(v',w)$, as required. 
\end{proof}

\sset\subsubsection{}\label{subsec_special-investig}
We wish next to investigate the specializations in the
valuation spectrum. Thus, let $A$ be any ring,
$v:A\to\Gamma_{\!\circ}$ a valuation. The convex subgroup
of $\Gamma_{\!v}$ generated by
$\Img\,(v)\setminus\Gamma^+_{\!\circ}$ is called the
{\em characteristic subgroup} of $v$, and shall
be denoted
$$
c\Gamma_{\!v}.
$$
Now, let $\Delta\subset\Gamma$ be any convex subgroup,
and $\pi:\Gamma\to\Gamma/\Delta$ the projection; we
associate with $\Delta$ two mappings
$$
\begin{aligned}
v_\Delta :\, & A\to(\Gamma/\Delta)_\circ
& & \qquad
a\mapsto\pi_\circ\circ v(a) \\
v^\Delta :\, & A\to\Delta_\circ
& & \qquad
a\mapsto\left\{\begin{array}{ll}
               v(a) & \text{if $v(a)\in\Delta$} \\
                0  & \text{otherwise}.
               \end{array}\right.
\end{aligned}
$$
\begin{lemma}\label{lem_specialize-vals}
With the notation of \eqref{subsec_special-investig},
the following holds :
\begin{enumerate}
\item
$v_\Delta$ is a valuation on $A$.
\item
$v^\Delta$ is a valuation if and only if
$c\Gamma_{\!v}\subset\Delta$.
\end{enumerate}
\end{lemma}
\begin{proof}(i) is clear. For (ii), suppose first that
$c\Gamma_{\!v}\subset\Delta$, pick any $a,b\in A$, and let us
show that $v(a+b)^\Delta\leq\max(v(a)^\Delta,v(b)^\Delta)$. This
is clear if both $v(a),v(b)\in\Delta_\circ$, so we may assume
that $v(b)\notin\Delta_\circ$, in which case notice that $v(b)<1$,
since $c\Gamma_{\!v}\subset\Delta$. Now, if $v(a)\in\Delta$,
then we must have $v(b)<v(a)$, as $\Delta$ is convex; therefore,
$v(a+b)=v(a)$ (remark \ref{rem_semi-norm}(iii)), and the
contention follows easily in this case. Lastly, if neither
of $v(a),v(b)$ lies in $\Delta$, the foregoing yields
$v(a),v(b)<1$, hence $v(a+b)<1$, and moreover we may assume
that $v(a+b)\leq v(a)$; then $v(a+b)$ cannot lie in $\Delta$,
as the latter is convex, so $v(a+b)^\Delta=0$, and the
assertion follows also in this case. Next, let us check
that $v(ab)^\Delta=v(a)^\Delta\cdot v(b)^\Delta$. Again, this
is clear if both $v(a),v(b)\in\Delta_\circ$, hence suppose
that $v(b)\notin\Delta_\circ$; if $v(a)\in\Delta_\circ$, then
it follows easily that $v(ab)\notin\Delta$, in which case
the assertion holds. If neither of $v(a),v(b)$ is in
$\Delta_\circ$, we have already remarked that $v(a),v(b)<1$,
hence $v(ab)<v(a)<1$, so $v(ab)$ cannot lie in $\Delta$,
as the latter is convex. Thus, the assertion holds also
in this case.

Conversely, suppose that $c\Gamma_{\!v}\not\subset\Delta$;
then there exists $a\in A$ such that $v(a)>1$ but
$v(a)\notin\Delta$. Set $b:=1-a$; it follows that
$v(b)=v(a)$, and $v(a+b)=1\in\Delta$, so $v(a+b)^\Delta=1$
but $v(a)^\Delta=v(b)^\Delta=0$, and therefore $v^\Delta$
is not a semi-norm.
\end{proof}

\begin{example} In the situation of example
\ref{ex_toric-valuations}, suppose that $A$ is a domain
and $P$ is integral with $P^\gp$ torsion-free, so that
$v_\phi$ is a valuation on $A[P]$, and let
$\Delta\subset\Gamma$ be a convex subgroup. 

(i)\ \
Lemma \ref{lem_specialize-vals}(ii) implies that
$v^\Delta_\phi$ is a valuation on $A[P]$ if and only if
$\phi(P)\subset\Gamma^+_\circ\cup\Delta$. Moreover, if
the latter condition holds, the mapping
$$
\phi^\Delta:P\to\Gamma_\circ
\qquad
x\mapsto\left\{\begin{array}{ll}
               \phi(x) & \text{if $\phi(x)\in\Delta$} \\
                0  & \text{otherwise}.
               \end{array}\right.
$$
is a morphism of monoids, and $v^\Delta_\phi=v_{\phi^\Delta}$.

(ii)\ \
Furthermore, $(v_\phi)_\Delta=v_{\pi\circ\phi}$, where
$\pi:\Gamma\to\Gamma/\Delta$ is the projection.
\end{example}

\begin{remark}\label{rem_special-investig}
With the notation of \eqref{subsec_special-investig},
we notice :

(i)\ \
Let $\bar v:\kappa(v)\to\Gamma_{\!\circ}$ be the residual
valuation of $v$, and $\kappa(v)^+$ the valuation ring
of $\bar v$; the correspondence of remark
\ref{rem_valuations}(vii) assigns to $\Delta$ a prime
ideal $\fp(\Delta)\subset\kappa(v)^+$, and we have
$$
\Ker\,v_\Delta=\Ker\,v
\qquad
\kappa(v_\Delta)^+=\kappa(v)^+_{\fp(\Delta)}=
\bar v{}^{-1}(\Gamma^+\cdot\Delta_\circ).
$$

(ii)\ \
Let $\pi:A\to\kappa(v)$ be the natural map; then
$c\Gamma_{\!v}\subset\Delta$ if and only if
$\pi(A)\subset\kappa(v)^+_{\fp(\Delta)}$.

(iii)\ \
Suppose now that $c\Gamma_{\!v}\subset\Delta$; then :
$$
\kappa(v^\Delta)\subset\Frac\,\kappa(v)^+/\fp(\Delta)
\qquad\text{and}\qquad
\kappa(v^\Delta)^+=\kappa(v^\Delta)\cap\kappa(v)^+/\fp(\Delta).
$$

(iv)\ \
The valuation $v_\Delta$ is a generization of $v$ in
$\Spv\,A$. Indeed, if $a,b\in A$ are any two
elements such that $v(a)\leq v(b)$, then clearly
$v_\Delta(a)\leq v_\Delta(b)$, whence the claim.

(v)\ \
Also, if $c\Gamma_{\!v}\subset\Delta$, then $v^\Delta$
is a specialization of $v$ in $\Spv\,A$ : indeed, if
$v^\Delta(a)\leq v^\Delta(b)$, it is easily seen that
$v(a)\leq v(b)$.
\end{remark}

\begin{definition} Let $A$ be a ring, and $v,w\in\Spv\,A$
any two elements.
\begin{enumerate}
\item
We say that $w$ is a {\em primary specialization} of $v$
in $\Spv\,A$, if there exists a convex subgroup
$\Delta\subset\Gamma_{\!v}$ containing $c\Gamma_{\!v}$, and
such that $w$ is (equivalent to) $v^\Delta$. In this case,
we also say that $v$ is a {\em primary generization} of $w$.
\item
We say that $w$ is a {\em secondary generization} of $v$
if there exists a convex subgroup $\Delta\subset\Gamma_{\!v}$
such that $w$ is (equivalent to) $v_\Delta$. In this case, we
also say that $v$ is a {\em secondary specialization} of $w$.
\item
We say that $w$ is a {\em generalized primary specialization}
of $v$, if $w$ is either a primary specialization of $v$,
or else $c\Gamma_{\!v}=\{1\}$ and $w$ is a trivial ({\em i.e.}
rank zero) valuation of $A$ with
$\Ker\,v^{\{1\}}\subset\Ker\,w$. In this case, $w$ is a
specialization of $v^{\{1\}}$ : see remark \ref{rem_Spv-of-ring}(iii).
\end{enumerate}
\end{definition}

Recall that specializations of points define a partial
ordering on $\Spv\,A$ : see remark \ref{rem_specialize}(ii).
Let $v$ be any semi-norm on $A$, and $I\subset A$ any ideal;
for the following lemma, we shall say that $I$ is
{\em $v$-convex}, if the following holds.
For every $a\in A$ such that there exists $b\in I$ with
$v(a)\leq v(b)$, we have $a\in I$ as well.

\begin{lemma}\label{lem_Spezializerungen}
Let $A$ be any ring, $v,w,u\in\Spv\,A$ any three elements.
We have :
\begin{enumerate}
\item
If $w$ is a primary specialization of $v$ and $u$ is a
primary specialization of $w$, then $u$ is a primary
specialization of $v$.
\item
If $w$ is a secondary generization of $v$ and $u$ is a
secondary generization of $w$, then $u$ is a secondary
generization of $v$.
\item
Suppose that $w$ is a generization of $v$. Then $w$ is
a secondary generization of $v$ if and only if\/
$\Ker\,w=\Ker\,v$.
\item
The primary specializations (resp. the secondary
generizations) of $v$ form a totally ordered subset
of the partially ordered set\/ $\Spv\,A$.
\item
The supports of the primary specializations of $v$ are
the $v$-convex prime ideals of $A$.
\end{enumerate}
\end{lemma}
\begin{proof}(i) and (ii) are clear from the explicit
description in \eqref{subsec_special-investig}.

(iii): Suppose that $w$ is a generization of $v$ whose
support equals that of $v$, and let $\bar v$ and $\bar w$
be the respective residual valuations on
$\kappa:=\kappa(v)=\kappa(w)$; from \eqref{eq_fibers-of-supp}
we deduce that $\bar w$ is a generization of $\bar v$
in $\Spv\,\kappa$, which means that
$\kappa(\bar v)^+\subset\kappa(\bar w)^+$. Then, let
$\Delta\subset\Gamma_{\!v}$ be the image of
$(\kappa(w)^+)^\times$; it follows easily that $\Delta$
is a convex subgroup, and taking into account remark
\ref{rem_valuations}(iv), we see that $\Gamma_{\!w}$ is
naturally identified with $\Gamma_{\!v}/\Delta$, and
$w=v_\Delta$.

(iv): Indeed, the explicit description of
\eqref{subsec_special-investig} implies more precisely
that the partially ordered set of secondary generizations
of $v$ is naturally identified with $\Spec\,\Gamma_{\!v}$,
where the latter is totally ordered by inclusion (see
definition \ref{def_ordered-group}(ii)). Likewise, we have
a surjective order-preserving map from
$\Spec\,\Gamma_{\!v}/c\Gamma_{\!v}$ onto the set of primary
specializations of $v$.

(v): Let $\Delta\subset\Gamma_{\!v}$ be a convex subgroup
containing $c\Gamma_{\!v}$; let also $a\in A$ and
$b\in\Ker(v^\Delta)$ with $v(a)\leq v(b)$. Then
$v(b)\notin\Delta$, and in particular $v(b)<1$, since
$c\Gamma_{\!v}\subset\Delta$; since $\Delta$ is convex, it
follows that $v(a)\notin\Delta$, hence $a\in\Ker(v^\Delta)$.
This shows that $\Ker(v^\Delta)$ is $v$-convex. Conversely,
suppose that the prime ideal $\fp\subset A$ is $v$-convex.
Hence, $v(x)\neq 0$ for every $x\in A\setminus\fp$, and
we let $\Delta\subset\Gamma_{\!v}$ be the subgroup generated
by $v(A\setminus\fp)$. We claim that $\gamma>v(x)$ for
every $\gamma\in\Delta$ and every $x\in\fp$. Indeed,
suppose that this condition fails for some $\gamma\in\Delta$
and $x\in\fp$; since $A\setminus\fp$ is a submonoid of
$(A,\cdot)$, we may write $\gamma=v(a)\cdot v(b)^{-1}$ for
some $a,b\in A\setminus\fp$, and we get $v(a)\leq v(bx)$,
which is absurd, since $\fp$ is $v$-convex. Now, let
$\Delta'\subset\Gamma_{\!v}$ be the smallest convex subgroup
containing $\Delta$; it is easily seen that we still
have $\gamma>v(x)$ for every $\gamma\in\Delta'$ and
every $x\in\fp$. On the other hand, by construction
we have $v(A\setminus\fp)\subset\Delta'$, so
$\fp=\{x\in A~|~\gamma>v(x)\ \text{for every $\gamma\in\Delta'$}\}$.
Moreover, say that $a\in A$ and $v(a)\geq 1$; then
$a\in A\setminus\fp$, since $\fp$ is $v$-convex, and
consequently $v(a)\in\Delta'$, {\em i.e.} the characteristic
subgroup of $v$ lies in $\Delta'$, so that
$v^{\Delta'}\in\Spv\,A$, and to conclude, it suffices to
remark that $\Ker\,v^{\Delta'}=\fp$.
\end{proof}

\begin{proposition}\label{prop_prprties-of-sigma}
Let $A$ be any ring, $v$ a valuation of $A$, and $\fp\subset A$
a prime ideal. We have :
\begin{enumerate}
\item
If $\fp\subset\Ker\,v$, there exists a primary generization
$w$ of $v$ with $\Ker\,w=\fp$. Especially, the map
$\sigma_A$ of remark {\em\ref{rem_Spv-of-ring}(iii)} is generizing.
\item
If $v\in\Spv^+A$, then every primary specialization and every
primary generization of $v$ in $\Spv\,A$ lies also in $\Spv^+A$,
and if $w$ is any primary specialization (or generization) of $v$,
then $\sigma^+_A(w)=\sigma^+_A(v)$ (notation of remark
{\em\ref{rem_Spv-of-ring}(v)}).
\item
If $v\in\Spv^+A$ and $\fp$ contains the center of $v$, there
exists a secondary specialization $w\in\Spv^+A$ of $v$ with
center equal to $\fp$. Especially, the map $\sigma^+_A$ is
specializing.
\end{enumerate}
\end{proposition}
\begin{proof}(i): The valuation $v$ factors through the projection
$A\to A':=A/\fp$ and a valuation $v'$ of $A'$; after replacing $A$
and $v$ by $A'$ and $v'$, we may then assume that $A$ is a domain
and $\fp=0$. In this case, let $\fq:=\Ker\,v\subset A$; according
to theorem \ref{th_val-cornerstone} there exists a valuation ring
$(V,\fm_V)$ of $K:=\Frac\,A$ such that $A\subset V$ and
$\fq=\fm_V\cap A$; by the same token, there exists a valuation
ring $(\bar W,\fm_{\bar W})$ of $V/\fm_V$ such that
$\bar W\cap\kappa(\fq)=\kappa(v)^+$. Let $\pi_V:V\to V/\fm_V$
be the projection; then $W:=\pi_V^{-1}(\bar W)$ is a valuation
ring of $K$; let $w_K$ be the unique valuation of $K$ with
$\kappa(w_K)^+=W$, and $w$ the restriction of $w_K$ to the
subring $A$. Notice that $\fm_V=\Ker(W\to\bar W)$ is a prime
ideal of $W$; then the bijection of remark \ref{rem_valuations}(vii)
assigns to $\fm_V$ a convex subgroup $\Delta(\fm_V)$ such that
$\fm_V=\{x\in W\ |\ \gamma<w_K(x)\
\text{for every}\ \gamma\in\Delta(\fm_V)\}$. Especially,
we see that $\fm_W$ is a $w_K$-convex ideal, hence $\fq$
is $w$-convex. Let $\Delta\subset\Gamma_{\!w}$ be the convex
subgroup generated by $w(A\setminus\fq)$; the proof of lemma
\ref{lem_Spezializerungen}(v) shows that
$c\Gamma_{\!w}\subset\Delta$ and $\Ker(w^\Delta)=\fq$. We need
to check that $w^\Delta$ is equivalent to $v$. Notice first
that for every $a,b\in A$ we have :
\set\begin{equation}\label{eq_je-trace}
w(a)\leq w(b)\Leftrightarrow a\in bW
\qquad\text{and}\qquad
v(a)\leq v(b)\Leftrightarrow\pi_V(a)\in\pi_V(b)\bar W.
\end{equation}
Now, suppose that $w^\Delta(a)\leq w^\Delta(b)$, and let us show
that $v(a)\leq v(b)$. Indeed, if $w^\Delta(b)\neq 0$, we have
$w^\Delta\in R_A\bigl(\frac{a}{b}\bigr)$, whence
$w\in R_A\bigl(\frac{a}{b}\bigr)$, so $w(a)\leq w(b)$,
and thus $v(a)\leq v(b)$, in light of \eqref{eq_je-trace}.
If $w^\Delta(b)=0$, we have as well $w^\Delta(a)=0$, and then
$v(a)=v(b)=0$, since $\Ker(v)=\Ker(w^\Delta)$. Conversely,
suppose that $v(a)\leq v(b)$; if $v(b)=0$, we see as in
the foregoing that $w^\Delta(a)=w^\Delta(b)=0$. Thus, let
$v(b)\neq 0$, and suppose by contradiction that
$w^\Delta(a)>w^\Delta(b)$; we deduce that
$w^\Delta(a),w^\Delta(b)\neq 0$, hence $a,b\neq 0$ and
$w(a)<w(b)$. Let $\fm_W=\pi_V^{-1}(\fm_{\bar W})$ be the
maximal ideal of $W$; we get $b\in a\cdot\fm_W$, so that
$\pi_V(b)\in\pi_V(a)\cdot\fm_{\bar W}$, {\em i.e.}
$v(a)>v(b)$, against our assumption. We may now apply
remark \ref{rem_valuations}(vi) to conclude.

(ii) is clear from the definitions.

(iii): Let $V\subset\kappa(v)$ be the valuation ring of
the residual valuation of $v$, and $\fm_V$ the maximal
ideal of $V$; then $v$ induces a ring homomorphism
$\pi:A\to V$ such that $\fq:=\pi^{-1}\fm_V$ is the center
of $v$. We deduce an injective ring homomorphism
$\bar A:=A/\fq\to\kappa:=V/\fm_V$, and we denote by
$\bar\fm\subset\bar A$ the image of $\fm$. By
corollary \ref{cor_cornerstone}, there exists a
valuation ring $\bar W$ of $\kappa$ that dominates the
image of $\bar A_{\bar\fm}$, and we denote by $W\subset V$
the preimage of $\bar W$. Then $W$ is a valuation ring
of $\kappa(v)$, and it corresponds to a secondary
specialization of $v$ in $\Spv^+A$.
\end{proof}

\begin{proposition}\label{prop_spezialiserungen}
Let $A$ be any ring, and $v\in\Spv\,A$ any element.
Then, every specialization of $v$ is a secondary
specialization of a generalized primary specialization
of $v$.
\end{proposition}
\begin{proof} Let $w$ be a specialization of $v$. If
$c\Gamma_{\!v}=\{1\}$ and $v(a)\geq 1$ for every
$a\in A\setminus\Ker\,w$, then $\Ker\,v^{\{1\}}\subset\Ker\,w$,
so the trivial semi-norm $u$ with support $\Ker\,w$ is
a generalized primary specialization of $v$, and $w$
is a secondary specialization of $u$. We may therefore
assume that either $c\Gamma_{\!v}\neq\{1\}$, or
$v(a)<1$ for some $a\in A\setminus\Ker\,w$.

\begin{claim}\label{cl_beobachtung}
Let $a,b\in A$ be any two elements with $v(a)\leq v(b)$,
$w(a)\neq 0$ and $w(b)=0$. Then $v(a)=v(b)\neq 0$. 
\end{claim}
\begin{pfclaim} The assumptions imply that
$w\in R_A(\frac{b}{a})$, and since $v$ is a generization
of $w$, it follows that $v\in R_A(\frac{b}{a})$ as well,
whence the claim.
\end{pfclaim}

\begin{claim}\label{cl_v-convex}
$\Ker\,w$ is $v$-convex.
\end{claim}
\begin{pfclaim} Indeed, let $x,y\in A$ be two elements with
\set\begin{equation}\label{eq_beobachtung}
v(x)\leq v(y)
\qquad\text{and}\qquad
w(y)=0
\end{equation}
and suppose, by way of contradiction, that $w(x)\neq 0$;
by claim \ref{cl_beobachtung} it follows that
\set\begin{equation}\label{eq_begrundung}
v(x)=v(y)\neq 0.
\end{equation}
Consider first the case where $c\Gamma_{\!v}\neq\{1\}$, and pick
$a\in A$ with $v(a)>1$; combining with \eqref{eq_beobachtung}
we get $v(x)\leq v(ay)$ and $w(ay)=0$, so claim
\ref{cl_beobachtung} yields $v(ay)=v(x)$, which contradicts
\eqref{eq_begrundung}. Next, if $v(a)<1$ for some $a\in A$
such that $w(a)\neq 0$, we get $v(ax)\leq v(y)$ and
$w(ax)\neq 0$, again by combining with \eqref{eq_beobachtung},
and claim \ref{cl_beobachtung} yields $v(ax)=v(y)$, again
a contradiction.
\end{pfclaim}

By claim \ref{cl_v-convex} and lemma \ref{lem_Spezializerungen}(iv,v),
there exists a unique primary specialization $u$ of $v$ with
$\Ker\,u=\Ker\,w$, and in view of lemma \ref{lem_Spezializerungen}(iii),
it suffices to show that $w$ is a specialization of $u$.
To this aim, let $a,b\in A$ be any two elements such that
$w\in R_A(\frac{a}{b})$; since $v$ is a generization of
$w$, we have $v\in R_A(\frac{a}{b})$ as well, and since $u$
is a primary specialization of $v$, we deduce $u(a)\leq u(b)$.
Lastly, since $\Ker\,u=\Ker\,w$ and $w(b)\neq 0$, we get
$u(b)\neq 0$, so $u\in R_A(\frac{a}{b})$, and the assertion
follows.
\end{proof}

\begin{lemma}\label{lem_sorite-Spezialiserung}
Let $A$ be a ring, $v,w\in\Spv\,A$ be two elements, such
that $w$ is a primary specialization of $v$. We have :
\begin{enumerate}
\item
For every secondary specialization $v'$ of $v$ there exists
a unique secondary specialization $w'$ of $w$ such that
$w'$ is a primary specialization of $v'$.
\item
For every secondary specialization $w'$ of $w$ there
exists a secondary specialization $v'$ of $v$ such that
$w'$ is a primary specialization of $v'$.
\item
For every secondary generization $v'$ of $v$ there exists
a unique secondary generization $w'$ of $w$ such that
$w'$ is a generalized primary specialization of $v'$.
\item
For every secondary generization $w'$ of $w$ there exists
a secondary generization $v'$ of $v$ such that $w'$ is
a primary specialization of $v'$.
\end{enumerate}
\end{lemma}
\begin{proof}(i): The uniqueness of $w'$ follows from
lemma \ref{lem_Spezializerungen}(iii,iv). For the existence,
let $\Delta'$ be a convex subgroup of $\Gamma_{\!v'}$ such
that $v=v'_{\Delta'}$, and $\Sigma$ a convex subgroup of
$\Gamma_{\!v}=\Gamma_{\!v'}/\Delta'$ containing $c\Gamma_{\!v}$
and such that $w=v^\Sigma$.
Denote by $\Sigma'$ the unique convex subgroup of
$\Gamma_{\!v'}$ such that $\Delta'\subset\Sigma'$ and
$\Sigma'/\Delta'=\Sigma$; then $c\Gamma_{\!v'}\subset\Sigma'$
and we may take $w':=(v')^{\Sigma'}$.

(ii): By remark \ref{rem_special-investig}(ii), there
exists a prime ideal $\fp\subset\kappa(v)^+$ such that
$\kappa(w)\subset K:=\Frac\,\kappa(v)^+/\fp$ and
$\kappa(w)^+=\kappa(w)\cap K^+$, where $K^+:=\kappa(v)^+/\fp$.
On the other hand, $\kappa(w')=\kappa(w)$, and
$\kappa(w')^+\subset\kappa(w)^+$. Let $\fm_K$, $\fm_w$
and $\fm_{w'}$ be the maximal ideals of respectively
$K^+$, $\kappa(w)^+$ and $\kappa(w')^+$; it follows
easily that $\fm_K\cap\kappa(w)=\fm_w\subset\fm_{w'}$.
Then $\fm_w$ is an ideal of $\kappa(w')^+$, and we have
injective ring homomorphisms
$$
B:=\kappa(w')^+/\fm_w\to\kappa(w)^+/\fm_w\to K^+/\fm_K.
$$
According to corollary \ref{cor_cornerstone}, we may
then find a valuation ring $\bar V$ of $K^+/\fm_K$ that
dominates $B$.
Since $B$ is a valuation ring, we easily deduce that
$B=\bar V\cap\Frac\,B$. Let $V\subset K^+$ be the preimage
of $\bar V$; then $V$ is a valuation ring of $K$ such that
$\kappa(w')^+=V\cap\kappa(w')$ (details left to the reader).
Let $V'\subset\kappa(v)^+$ be the preimage of $V$; then
$V'$ is a valuation ring of $\kappa(v)$, and we let $v'$
be the unique valuation of $A$ with $\Ker\,v'=\Ker\,v$ and
$\kappa(v')^+=V'$; it is easily seen that $v'$ is the sought
secondary specialization of $v$.

(iii): The uniqueness of $w'$ follows from
lemma \ref{lem_Spezializerungen}(iii,iv), and the
existence follows from proposition \ref{prop_spezialiserungen}.

(iv): Let $\Delta$ be a convex subgroup of $\Gamma_{\!v}$
containing $c\Gamma_{\!v}$ and such that $w=v^\Delta$, and
$\Sigma$ a convex subgroup of $\Delta$ such that $w'=w_\Sigma$.
Set $v':=v_\Sigma$; then $\Delta':=\Delta/\Sigma$ is a convex
subgroup of $\Gamma_{\!v'}=\Gamma_{\!v}/\Sigma$ that contains
$c\Gamma_{\!v'}$, and $w'=(v')^{\Delta'}$.
\end{proof}

\begin{proposition}\label{prop_reverse-order}
Let $A$ be a ring, and $v\in\Spv\,A$ any element. Then
every specialization of $v$ is a primary specialization
of a secondary specialization of $v$.
\end{proposition}
\begin{proof} Let $w$ be a specialization of $v$. If $w$
is a secondary specialization of a primary specialization
of $v$, the assertion follows from lemma
\ref{lem_sorite-Spezialiserung}(ii). In light of proposition
\ref{prop_spezialiserungen}, we may then assume that the
characteristic subgroup of $v$ is trivial and
$\Ker\,v^{\{1\}}\subset\Ker\,w$. In this case, proposition
\ref{prop_prprties-of-sigma}(i) yields a primary generization
$w'$ of $w$ with $\Ker\,w'=\Ker\,v^{\{1\}}$. Clearly
$w'$ is a secondary specialization of $v^{\{1\}}$, so
lemma \ref{lem_sorite-Spezialiserung}(ii) ensures that
we may find a secondary specialization $v'$ of $v$
which is also a primary generization of $w'$, and
then $w$ is a primary specialization of $v'$ (lemma
\ref{lem_Spezializerungen}(i)).
\end{proof}

\begin{lemma}\label{lem_image-of-special}
Let $f:A\to B$ be any ring homomorphism, and
$v,w\in\Spv\,B$ any two elements. We have :
\begin{enumerate}
\item
If $w$ is a primary (resp. secondary) specialization of
$v$, then $\Spv(f)(w)$ is a primary (resp. secondary)
specialization of\/ $\Spv(f)(v)$.
\item
If\/ $\Spec\,f$ is an open immersion, and $\Spv(f)(w)$ is
a primary (resp. secondary) specialization of\/ $\Spv(f)(v)$,
then $w$ is a primary (resp. secondary) specialization of $v$.
\item
$\Spv\,f$ restricts to a surjection from the set of
secondary specializations (resp. generizations) of $v$ onto
the set of secondary specializations (resp. generizations)
of\/ $\Spv(f)(v)$.
\item
Suppose that $w$ is a primary specialization of $v$, and
denote by $P$ (resp. $Q$) the set of primary specializations
of $v$ in $\Spv\,B$ (resp. of\/ $\Spv(f)(v)$ in $\Spv\,A$)
that are also primary generizations of $w$ in $\Spv\,B$
(resp. of\/ $\Spv(f)(w)$ in $\Spv\,A$). Then $\Spv\,f$
restricts to a surjection $P\to Q$.
\end{enumerate}
\end{lemma}
\begin{proof}(i) is immediate from the definitions.

(ii): Set $v':=\Spv(f)(v)$; if $\Spec\,f$ is an open
immersion, $f$ induces an isomorphism
$\kappa(v')\isom\kappa(v)$, and similarly for $w$. The
assertion is an immediate consequence.

(iii): Clearly, $f$ induces an injective map
$i:\Gamma_{\!v'}\to\Gamma_{\!v}$, and $\Spec\,i$ is surjective,
by remark \ref{rem_ordered-gps}(v), so $\Spv\,f$ maps
surjectively the secondary generizations of $v$ onto the
ones of $v'$. Next, let $w'$ be a secondary specialization
of $v'$ in $\Spv\,A$; then $\kappa(w')=\kappa(v')$ and
$\kappa(w')^+$ is a valuation ring of $\kappa(v')$ contained
in $\kappa(v')^+$. Let $\fm\subset\kappa(v)^+$ be the maximal
ideal, so that $\fm':=\fm\cap\kappa(v')$ is the maximal ideal
of $\kappa(v')^+$ and a prime ideal of $\kappa(w')^+$; then
$K^+:=\kappa(w')^+/\fm'$ is a valuation ring of
$K:=\kappa(v')^+/\fm'$, and by corollary \ref{cor_cornerstone}
there exists a valuation ring $V$ of $\kappa(v)^+/\fm$ that
dominates $K^+$. The preimage of $V$ in $\kappa(v)^+$
is a valuation ring of $\kappa(v)$ and the corresponding
valuation $w$ of $B$ is a secondary specialization of
$v$ with $\Spv(f)(w)=w'$.

(iv): Again, since the induced morphism of ordered groups
$i$ is injective, the assertion follows easily from remark
\ref{rem_ordered-gps}(v).
\end{proof}

\begin{theorem} Let $f:A\to B$ be a flat ring homomorphism,
and $v\in\Spv\,B$ any element. Then the following holds :
\begin{enumerate}
\item
$\Spv\,f$ restricts to a surjection from the set of primary
generizations of $v$ in $\Spv\,B$ to the set of primary
generizations of\/ $\Spv(f)(v)$ in $\Spv\,A$.
\item
$\Spv\,f$ is generizing.
\end{enumerate}
\end{theorem}
\begin{proof}(i): Let $t$ be any primary generization of
$s:=\Spv(f)(v)$. According to remark
\ref{rem_special-investig}(i), there exists a prime ideal
$\fp\subset\kappa(t)^+$ such that
\begin{itemize}
\item
the image of the natural map $A\to\kappa(t)$ lies
in $A':=\kappa(t)^+_\fp$
\item
$\kappa(s)$ is naturally identified with a subfield
of the residue field $\bar\kappa:=A'/\fp A'$ of $A'$
\item
under this identification, we have
$\kappa(s)^+=\kappa(s)\cap(\kappa(t)^+/\fp)$.
\end{itemize}
With this notation, let $h:A\to A'$ the resulting map, and
$$
\bar t:\kappa(t)=\Frac\,A'\to\Gamma_{\!t\circ}
\qquad\text{and}\qquad
\bar s:\bar\kappa\to\Gamma_{\!s\circ}
$$
be respectively the residual valuation of $t$ and the
valuation of $\bar\kappa$ such that
$\kappa(\bar s)^+=\kappa(t)^+/\fp$. Then, $\bar t$ and
$\bar s$ yield elements $t'$ and respectively $s'$ of
$\Spv\,A'$ such that
$$
s=\Spv(h)(s')
\qquad\text{and}\qquad
t=\Spv(h)(t').
$$
Let $f':A'\to B':=A'\otimes_AB$ be the induced map; by
proposition \ref{prop_tensor-Spv} we may find $v'\in\Spv\,B'$
whose image in $\Spv\,A'$ (resp. in $\Spv\,B$) equals $s'$
(resp. $v$), and in light of lemma \ref{lem_image-of-special}(i),
it suffices to exhibit a primary generization $w'$ of $v'$
in $\Spv\,B'$ with $\Spv(f')(w')=t'$.
However, since $f'$ is flat, there exists a prime ideal
$\fr\subset\fq:=\Ker\,v'$ in $B'$, such that $f'^{-1}\fr=0$
(\cite[Th.9.5]{Mat}); we set $B'':=B'_\fq/\fr B'_\fq$ and
we let $f'':A'\to B''$ be the resulting map.
There exists a unique valuation $v''$ of $B''$
whose image in $\Spv\,B'$ equals $v'$, and it suffices
to exhibit a primary generization $w''$ of $v''$ in
$\Spv\,B''$ such that $\Spv(f'')(w'')=t'$. To this aim,
set $K'':=\Frac\,B''$; by corollary \ref{cor_cornerstone}
we may find a valuation ring $V$ of $K''$ that dominates
$B''$. Let $\kappa''$ be the residue field of $V$, and
denote by $R$ the image of $\kappa(v'')^+$ in $\kappa''$;
we may find a valuation ring $V'$ of $\kappa''$ that
dominates $R$ (corollary \ref{cor_cornerstone}). Denote
by $V''\subset V$ the preimage of $V'$. Then $V''$ is a
valuation ring of $K''$, and we may take for $w''$ the
unique valuation of $B''$ with $\kappa(w'')^+=V''$.

(ii) is an immediate consequence of (i), lemma
\ref{lem_image-of-special}(iii) and proposition
\ref{prop_reverse-order}.
\end{proof}

\begin{theorem} Let $f:A\to B$ be an integral ring
homomorphism, $v\in\Spv\,B$ any element. The following
holds :
\begin{enumerate}
\item
If\/ $v'\in\Spv\,B$ is a specialization of $v$ with
$\Spv(f)(v')=\Spv(f)(v)$, then $v=v'$.
\item
$\Spv\,f$ restricts to a bijection from the set of primary
specializations of $v$ in $\Spv\,B$ to the set of primary
specializations of\/ $\Spv(f)(v)$ in $\Spv\,A$.
\item
$\Spv\,f$ is specializing.
\item
If $f$ is injective, $\Spv\,f$ is surjective.
\item
If $A$ is integral and normal, $B$ is integral and $f$
is injective, then $\Spv\,f$ restricts to a surjection
from the set of primary generizations of $v$ in $\Spv\,B$
onto the set of primary generizations of\/ $\Spv(f)(v)$ in
$\Spv\,A$.
\end{enumerate}
\end{theorem}
\begin{proof} Set $w:=\Spv(f)(v)$, and let $\Gamma_{\!v}$
and $\Gamma_{\!w}$ be the value groups of $v$ and $w$.

(i): Let $\fp$ and $\fq$ be the supports of $v$ and
respectively $v'$; under the current assumptions, we
have $f^{-1}\fp=f^{-1}\fq$, and $\fq$ is a specialization
of $\fp$ in $\Spec\,B$; then $\fq=\fq'$, by
\cite[Th.9.3(ii)]{Mat}. Taking into account proposition
\ref{prop_reverse-order}, we conclude that $v'$ is a
secondary specialization of $v$. However, since $f$ is
integral, it induces an isomorphism
$\Gamma_{\!v}\otimes_\Z\Q\isom\Gamma_{\!w}\otimes_\Z\Q$.
On the other hand, we have a natural bijection between
$\Spec\,\Gamma_{\!v}$ and the set $S_v$ of secondary
specializations of $v$ in $\Spv\,B$, and likewise for the
set $S_w$ of secondary specializations of $w$; taking into
account remark \ref{rem_ordered-gps}(v), we deduce that
$f$ induces a bijection $S_v\isom S_w$, whence the assertion.

(ii): We remark :

\begin{claim}\label{cl_same-char-subgroup}
Under the assumptions of the theorem, let $b\in B$
be any element. Then there exists $a\in A$ such that
$v(b)\leq w(a)$.
\end{claim}
\begin{pfclaim} Since $f$ is integral, we may find
$a_1,\dots,a_n\in A$ such that $b^n+a_1b^{n-1}+\cdots+a_n=0$.
It follows that $v(b)^n\leq v(b)^i\cdot w(a_{n-i})$ for
some $i<n$. If $v(b)=0$, there is nothing to prove;
otherwise, we get $v(b)^{n-i}\leq w(a_{n-i})$. If $v(b)\leq 1$,
there is nothing to prove; otherwise, we deduce that
$v(b)\leq w(a_{n-i})$, as required.
\end{pfclaim}

From claim \ref{cl_same-char-subgroup} it follows easily
that $c\Gamma_{\!w}=c\Gamma_{\!v}\cap\Gamma_{\!w}$, whence
the assertion.

(iii) follows immediately from (ii), proposition
\ref{prop_reverse-order} and lemma
\ref{lem_image-of-special}(iii).

(iv): Notice first that $\Spec\,f$ is surjective, under
these assumptions. Indeed, let $\fp\subset A$ be any prime
ideal, and pick a minimal prime ideal $\fq\subset A$
contained in $\fp$; the induced ring homomorphism
$f_\fq:A_\fq\to B_\fq$ is still injective, so $B\neq 0$;
now, if $\fq'\subset B_\fq$ is any prime ideal, then
$f_\fq^{-1}(\fq')=\fq$, and by \cite[Th.9.4(i)]{Mat} there
exists a specialization $\fp'$ of $\fq'$ in $\Spec\,B$ such
that $\Spec\,(f)(\fp')=\fp$. Now, let $w\in\Spv\,A$ be any
element, and $\fp'\in\Spec\,B$ any prime ideal such that
$f^{-1}\fp'=\fp:=\Ker\,w$; it follows easily from theorem
\ref{th_val-cornerstone} that there exists a valuation
$w'$ of the residue field $\kappa(\fp')$ whose restriction
to $\kappa(\fp)$ is equivalent to $v$.

(v): Let $w'$ be a primary generization of $w$ in $\Spv\,A$.
We need to exhibit a primary generization $v'$ of $v$ in
$\Spv\,B$ such that $\Spv(f)(v')=w'$. To this aim, pick
a normal algebraic field extension $E$ of $\Frac\,A$
containing $\Frac\,B$, and denote by $C$ the integral
closure of $A$ in $E$. By (iv), we may find $w'',v''\in\Spv\,C$
such that $w''_{|A}=w'$ and $v''_{|B}=v$. By (ii), we may
find a primary specialization $t$ of $w''$ in $\Spv\,C$
such that $t_{|A}=w$. By construction, $t_{|A}=v''_{|A}$,
so there exists an automorphism $h$ of the $A$-algebra
$C$ such that $\Spv(h)(t)=v''$
(\cite[Ch.V, \S2, no.3, Prop.6 and
Ch.VI, \S8, no.6, Cor.1]{BouAC}). It is easily
seen that $v':=\Spv(h)(w'')_{|B}$ will do.
\end{proof}

\begin{theorem}\label{th_Chevalley}
Let $f:A\to B$ be a finitely presented ring homomorphism,
and $T\subset\Spv\,B$ any constructible subset. Then
$\Spv(f)(T)$ is a constructible subset of\/ $\Spv\,A$.
\end{theorem}
\begin{proof} We may assume that there exist elements
$a_1,b_1,\dots,a_n,b_n,c_1,d_1,\dots,c_m,d_m\in B$ such
that $T$ is the set of all $v\in\Spv\,B$ with $v(a_i)<v(b_i)$
and $v(c_j)\leq v(d_j)$ for every $i=1,\dots,n$ and
$j=1,\dots,m$. We remark :

\begin{claim}\label{cl_ok-for-closed-immers}
The theorem holds if $\Spec\,f$ is a closed immersion.
\end{claim}
\begin{pfclaim} Indeed, in this case $f$ is a surjective
map, whose kernel is a finitely generated ideal $I\subset A$.
Pick a finite system $x_1,\dots,x_r\in A$ of generators of $I$,
and for every $i=1,\dots,n$, $j=1,\dots,m$ choose elements
$a'_i,b'_i,c'_j,d'_j\in A$ whose images equal respectively
$a_i,b_i,c_j,d_j$ in $B$. Then, $\Spv(f)(T)$ is the constructible
subset of $\Spv\,A$ of all valuations $v$ of $A$ such that
$v(x_k)=0$ for every $k=1,\dots,r$, and $v(a'_i)<v(b'_i)$
and $v(c'_j)\leq v(d'_j)$ for every $i=1,\dots,n$ and
$j=1,\dots,m$.
\end{pfclaim}

\begin{claim}\label{cl_Spv-plus}
(i)\ \
We may assume that $A$ and $B$ are $\Z$-algebras of
finite type.
\begin{enumerate}
\addenu
\item
For every $\Z$-algebra $R$ of finite type, $\Spv^+R$
is a rational subset of $\Spv\,R$ (notation of
remark \ref{rem_Spv-of-ring}(v)).
\end{enumerate}
\end{claim}
\begin{pfclaim}(i): Indeed, we may write $A$ as the union
of the filtered system $(A_\lambda~|~\lambda\in\Lambda)$
of its $\Z$-subalgebras of finite type, and then there
exists $\lambda\in\Lambda$ such that
$B=B_\lambda\otimes_{A_\lambda}A$ for a ring homomorphism
$f_\lambda:A_\lambda\to B_\lambda$ of finite type. Moreover,
after replacing $\lambda$ by a larger index, we may
assume that there exist elements
$a'_1,b'_1,\dots,a'_n,b'_n,c'_1,d'_1,\dots,c'_m,d'_m\in B_\lambda$
such that $a'_i\otimes 1=a_i$, $b_i=b'_i\otimes 1$,
$c_j=c'_j\otimes 1$ and $d_j=d'_j\otimes 1$ in $B$,
for every $i=1,\dots,n$ and $j=1,\dots,m$. We then
let $T'\subset\Spv\,B_\lambda$ be the constructible
subset of all $v\in\Spv\,B_\lambda$ such that $v(a'_i)<v(b'_i)$
and $v(c'_j)\leq v(d'_j)$ for every $i=1,\dots,n$ and
$j=1,\dots,m$. There follows a commutative diagram
of topological spaces
$$
\xymatrix{ \Spv\,B \ar[r]^-{\pi_B} \ar[d]_{\Spv\,f} &
\Spv\,B_\lambda \ar[d]^{\Spv\,f_\lambda} \\
\Spv\,A \ar[r]^-{\pi_A} & \Spv\,A_\lambda
}$$
clearly $T=\pi^{-1}_BT'$, and it follows easily from
proposition \ref{prop_tensor-Spv} that
$$
\pi_A^{-1}(\Spv(f_\lambda)(T'))=\Spv(f)(T).
$$
Since $\pi_A$ is spectral (remark \ref{rem_Spv-of-ring}(ii)),
we are then reduced to showing that the theorem holds
for the map $f_\lambda$.

(ii): Pick a finite system $x_1,\dots,x_r$ of elements
of $R$ such that $R=\Z[x_1,\dots,x_r]$; then it is easily
seen that $\Spv^+R=R_A(\frac{x_1}{1},\dots,\frac{x_r}{1})$.
\end{pfclaim}

Henceforth, we assume that $A$ and $B$ are $\Z$-algebras
of finite type, and we set
$$
C:=B[X_i,Y_j~|~i=1,\dots,n,\ j=1,\dots,m]/I
$$
where $I$ is the ideal generated by the system
$(a_i-X_ib_i,c_j-Y_jd_j~|~i=1,\dots,n,\ j=1,\dots,m)$.
We let $T'\subset\Spv\,C$ be the constructible subset
of all $w\in\Spv\,C$ such that
$$
0<w(b_i)
\qquad
w(X_i)<1
\qquad
w(Y_j)\leq 1
\qquad
\text{for $1\leq i\leq n$ and $1\leq j\leq m$}
$$
and denote by $g:B\to C$ the natural ring homomorphism.

\begin{claim}\label{cl_reduce-to-simple-ineq}
$\Spv(g)(T')=T$.
\end{claim}
\begin{pfclaim} Let $w\in T'$ be any element, and set
$v:=\Spv(g)(w)$; we may regard $w$ as a valuation of
$B[X_1,\dots,X_n,Y_1,\dots,Y_m]$ such that $w(b_i)>0$
and $w(a_i)=w(X_i)\cdot w(b_i)$ for every $i=1,\dots,n$,
and $w(c_j)=w(Y_j)\cdot w(d_j)$ for every $j=1,\dots,m$.
It follows easily that $v\in T$. Conversely, let $v\in T$
be any element, and set $D:=C\otimes_B\kappa(v)$.
Moreover, set $\Sigma:=\{1\leq j\leq m~|~v(d_j)=0\}$.
A simple inspection yields an isomorphism
$$
h:D\isom\kappa(v)[Y_j~|~j\in\Sigma]
$$
of $\kappa(v)$-algebras, such that :
\begin{itemize}
\item
$h(X_i)=h(a_i)\cdot h(b_i)^{-1}$ for $i=1,\dots,n$
\item
$h(Y_j)=Y_j$ for $j\in\Sigma$
\item
$h(Y_j)=h(c_j)\cdot h(d_j)^{-1}$ for every
$j\in\{1,\dots,m\}\setminus\Sigma$.
\end{itemize}
Let $(V,\fm_V)$ be any valuation ring of $\Frac\,D$
containing $\kappa(v)^+[Y_j~|~j\in\Sigma]$ and such
that $\fm_V\cap\kappa(v)$ is the maximal ideal of
$\kappa(v)^+$ (theorem \ref{th_val-cornerstone});
then $V$ corresponds to a valuation $\bar w$ of $D$
whose restriction to $\kappa(v)$ is equivalent to the
residual valuation $\bar v$ of $v$. So, $\bar w$
induces a valuation $w$ of $C$ with $\Spv(g)(w)=v$,
and by construction we have $w(X_i)=v(a_i)\cdot v(b_i)^{-1}<1$
for every $i=1,\dots,n$; likewise, we get
$w(Y_j)=v(c_j)\cdot v(d_j)^{-1}\leq 1$ for $j\in\Sigma$,
and $w(Y_j)\leq 1$ for $j\in\{1,\dots,m\}\setminus\Sigma$.
The claim follows.
\end{pfclaim}

In view of claim \ref{cl_reduce-to-simple-ineq}, we
may replace $B$ by $C$, and $T$ by $T'$, after which
we may assume:
\begin{itemize}
\item[(C)]
There exist finite sequences
$a_\bullet:=(a_1,\dots,a_n),b_\bullet:=(b_0,b_1,\dots,b_m)$ of
elements of $B$ such that
$$
T=R_B\Bigl(\frac{b_0}{b_0},\frac{b_1}{1},\dots,\frac{b_m}{1}\Bigr)
\setminus \Bigl(R_B\Bigl(\frac{1}{a_1}\Bigr)\cup\dots\cup
R_B\Bigl(\frac{1}{a_n}\Bigr)\Bigr).
$$
\end{itemize}
We shall argue by induction on the dimension of $B$.
Since $B$ is noetherian, $\Spec\,B$ has finitely many
irreducible components $Z_1,\dots,Z_k$, and
$$
\Spv\,B=\bigcup_{i=1}^k\Spv\,Z_i
$$
(notation of \eqref{subsec_Spv-of-scheme}); moreover,
since the inclusion map $Z_i\to\Spec\,B$ is spectral,
the subset $\Spv\,Z_i\cap T$ is constructible in
$\Spv\,Z_i$, for every $i=1,\dots,k$. Thus, it suffices
to show that the restriction $\Spv\,Z_i\to\Spv\,A$
of $\Spv\,f$ maps constructible subsets to constructible
subsets, for every $i=1,\dots,k$. We may therefore
assume from start that $\Spec\,B$ is also irreducible
and reduced, {\em i.e.} that $B$ is a domain. Then,
the topological closure $W$ of $\Spec(f)(\Spec\,B)$ in
$\Spec\,A$ is irreducible (lemma \ref{lem_irreducible}(v));
by virtue of claim \ref{cl_ok-for-closed-immers}, it
suffices to show that the induced map $\Spv\,B\to\Spv\,W$
sends constructible subsets to constructible subsets.
We may then assume as well that $A$ is a domain, and
$f$ is injective. In this situation, if $\dim\,B=0$,
we see that $A$ and $B$ are finite fields, and the
assertion is obvious. Hence, we suppose henceforth that
$d:=\dim\,B>0$, and that the theorem has already been shown
for every $\Z$-algebra $B$ of finite type of dimension $<d$.

Pick finite subsets $I\subset  A$ and $J\subset B$ such
that $A=\Z[I]$ and $B=A[J]$ and $I\cap J=\emptyset$. For
every pair of subsets $I'\subset I$ and $J'\subset J$ consider
the constructible open subsets
$$
\begin{aligned}
S_A(I'):= &\,
R_A\Bigl(\frac{x}{1}~|~x\in I'\Bigr)\cap
R_A\Bigl(\frac{1}{x}~|~x\in I\setminus I'\Bigr)\subset\Spv\,A \\
S_B(I',J'):= &\,
R_B\Bigl(\frac{x}{1}~|~x\in I'\cup J'\Bigr)\cap
R_B\Bigl(\frac{1}{x}~|~x\in(I\cup J)\setminus(I'\cup J')\Bigr)
\subset\Spv\,B.
\end{aligned}
$$
Clearly $\Spv\,A=\bigcup_{I'\subset I}S_A(I')$, and
likewise for $\Spv\,B$, so it suffices to show that
the subset $\Spv(f)(T\cap S_B(I',J'))$ is constructible
in $S_A(I')$, for any such $I'$ and $J'$ (lemma
\ref{lem_sorite-construct}(iii,x.a)). Now, for given $I'$
and $J'$, let :
\begin{itemize}
\item
$A_1:=\Z[I'\cup\{x^{-1}~|~x\in I\setminus I'\}]\subset\Frac\,A$
\item
$B_1:=
A_1[J'\cup\{x^{-1}~|~x\in J\setminus J'\}\cup a_\bullet\cup b_\bullet]
\subset\Frac\,B$
\item
$A_2:=A\cdot A_1\subset\Frac\,A$ and
$B_2:=B\cdot B_1\subset\Frac\,B$.
\end{itemize}
The natural inclusions maps $A\to A_2\leftarrow A_1$ and
$B\to B_2\leftarrow B_1$ are localizations, so $\dim B_1=\dim B_2=d$
(\cite[Ch.IV, Prop.10.6.1(ii)]{EGAIV-3}), and $f$ extends to
an injective ring homomorphism $f_2:A_2\to B_2$, which in turn
restricts to a map $f_1:A_1\to B_1$. So we get a commutative
diagram of continuous maps
$$
\xymatrix{
\Spv\,B_1 \ar[d]_{\Spv\,f_1} &
\ar[l] \Spv\,B_2 \ar[r] \ar[d]_{\Spv\,f_2} &
\Spv B \ar[d]^{\Spv\,f} \\
\Spv\,A_1 & \ar[l] \Spv\,A_2 \ar[r] & \Spv A
}$$
whose horizontal arrows are quasi-compact open immersions.
Moreover, $S_A(I')$ and $S_B(I',J')$ lie respectively in the
images of $\Spv\,A_2$ and $\Spv\,B_2$, and their images in
$\Spv\,A_1$ and $\Spv\,B_1$ are the subsets
$$
\begin{aligned}
Z_A(I'):= &\, \{v\in\Spv^+A_1~|~
\text{$v(x^{-1})\neq 0$ for every $x\in I\setminus I'$}\} \\
Z_B(I',J'):= &\, \{v\in\Spv^+B_1~|~
\text{$v(x^{-1})\neq 0$ for every $x\in(I\cup J)\setminus(I'\cup J')$}\}.
\end{aligned}
$$
Hence, we may replace $A$ by $A_1$, $B$ by $B_1$, and suppose that
$$
T=\Bigl(R_B\Bigl(\frac{b_0}{b_0}\Bigr)\setminus
R_B\Bigl(\frac{1}{b_1}\Bigr)\cup\cdots\cup
R_B\Bigl(\frac{1}{b_n}\Bigr)\Bigr)\cap\Spv^+B
\qquad
\text{for certains $b_0,\dots,b_n\in B$}.
$$
Next, by \cite[Partie I, Th.5.2.2]{Gr-Ra} we may find an
integer $r>0$ and elements $t_1,\dots,t_r\in A\setminus\{0\}$
such that, for every $i=1,\dots,r$ the inclusion map of
subrings of $B'_i:=B[t_i^{-1}]$
$$
A_i:=A\Bigl[\frac{t_1}{t_i},\dots,\frac{t_r}{t_i}\Bigr]\to
B_i:=B\Bigl[\frac{t_1}{t_i},\dots,\frac{t_r}{t_i}\Bigr]
$$
is flat. We consider the constructible subsets
$$
S_0:=\{v\in\Spv\,B~|~\text{$v(t_i)=0$ for $i=1,\dots,r$}\}
\qquad
S_i:=R_B\Bigl(\frac{t_1}{t_i},\dots,\frac{t_r}{t_i}\Bigr)
\qquad
\text{($i=1,\dots,r$)}.
$$
Since $\Spv\,B=\bigcup_{i=0}^rS_i$, it suffices to show
that $\Spv(f)(T\cap S_i)$ is constructible in $\Spv\,A$
for every $i=0,\dots,r$. However, $T\cap S_0$ is
a constructible subset of $\Spv\,B/I$, where $I\subset B$
is the ideal generated by $t_1,\dots,t_r$; since
$\dim\,B/I<\dim\,B$, the inductive assumption yields the
assertion for $i=0$. For $i=1,\dots,r$, set $A'_i:=A[t^{-1}_i]$;
we get a commutative diagram of continuous maps
$$
\xymatrix{
\Spv\,B'_i \ar[r]^-{\phi'_i} \ar[d] &
\Spv\,B_i \ar[r]^-{\phi_i} \ar[d] &
\Spv\,B \ar[d] \\
\Spv\,A'_i \ar[r]^-{\psi'_i} &
\Spv\,A_i \ar[r]^-{\psi_i} & \Spv\,A
}$$
and we notice that the maps $\phi'_i$, $\psi'_i$,
$\psi_i\circ\psi'_i$ and $\phi_i\circ\phi'_i$ are quasi-compact
open immersions, and $T\cap S_i=\phi_i(T_i)$, for the constructible
subsets of $\Spv\,B'_i$
$$
T_i:=
\Bigl(R_{B_i}\Bigl(\frac{t_ib_0}{t_ib_0}\Bigr)\setminus
\Bigl(R_{B_i}\Bigl(\frac{1}{b_1}\Bigr)\cup\cdots\cup
R_{B_i}\Bigl(\frac{1}{b_n}\Bigr)\Bigr)\Bigr)\cap\Spv^+B_i
\subset\phi'_i(\Spv\,B'_i).
$$
In light of lemma \ref{lem_sorite-construct}(iii,x.a), it then
suffices to show that the image of $T_i$ is constructible in
$\Spv\,A_i$. Thus, we may replace $A$, $B$, $b_0$ and $T$
respectively by $A_i$, $B_i$, $t_ib_0$ and $T_i$, and assume
furthermore that $f$ is a flat ring homomorphism.

Let $\phi:\Spec\,B/b_0B\to\Spec\,A$ be the restriction
of $\Spec\,f$; there exists a non-empty affine open
subset $U\subset\Spec\,A$ such that the restriction
$\phi^{-1}U\to U$ of $\phi$ is a flat morphism of schemes
(\cite[Ch.IV, Th.6.9.1]{EGAIV-2}). For every $\fp\in U$,
we have a short exact sequence of $A$-modules
$$
\cE\quad :\quad
0\to B_\fp\to B_\fp\to B_\fp/b_0B_\fp\to 0
$$
and the induced sequence $\cE\otimes_A\kappa(\fp)$ is
still exact, since $B_\fp/b_0B_\fp$ is a flat $A_\fp$-module
(notation of definition \ref{def_strict-loc}(i)). In
other words, we have :
\begin{itemize}
\item[(D)]
The image of $b_0$ in $B\otimes_A\kappa(\fp)$
is a regular element, for every $\fp\in U$.
\end{itemize}
Set
$$
Z:=\Spec\,A\setminus U
\qquad
Z_B:=\Spec(f)^{-1}Z
\qquad
U_B:=\Spec(f)^{-1}U.
$$
Then $\Spv\,Z$ is a constructible closed subset of $\Spv\,A$
and $Z_B$ is an affine scheme of dimension $<d$; also,
$T':=T\cap \Spv\,Z_B$ is a constructible subset of
$\Spv\,Z_B$, so our inductive assumption says that
$\Spv(f)(T')$ is constructible in $\Spv\,A$. It remains
therefore only to show that the subset
$\Spv(f)(T\cap\Spv\,U_B)$ is constructible in $\Spv\,A$.

Let $J\subset B$ be the ideal generated by $b_1,\dots,b_n$;
by \cite[Ch.IV, Th.1.8.4]{EGAIV}, the subset
$W:=\Spec(f)(\Spec\,B/J)$ is constructible in $\Spec\,A$,
hence
$$
\Spv^+W:=(\sigma^+_A)^{-1}(W)
$$
is constructible in $\Spv^+A$ (remark \ref{rem_Spv-of-ring}(v)).
In light of claim \ref{cl_Spv-plus}(ii), to conclude the proof,
it then suffices to show :

\begin{claim} Under the current assumptions, we have
$$
\Spv(f)(T\cap\Spv\,U_B)=\Spv^+W\cap\Spv\,U.
$$
\end{claim}
\begin{pfclaim}[] A simple inspection shows that
$\Spv(f)(T)\subset\Spv^+W$, so it remains only to check
that $\Spv^+W\cap\Spv\,U\subset\Spv(f)(T\cap\Spv\,U_B)$.
Now, let $v\in\Spv^+W\cap\Spv\,U$ be any element, and
let $\pi:A\to\kappa(v)$ be the induced map. Since
$v\in\Spv^+A$, we have $\Img\,\pi\subset V:=\kappa(v)^+$,
and we set $B':=B\otimes_AV$. Let also $\fm_V$ and $\eta_V$
be respectively the closed point and the generic point
of $\Spec\,V$; since $v\in\Spv\,U$, the image of $\eta_V$
in $\Spec\,A$ lies in $U$. By construction, there exists
$\fp\in W$ such that $\Spec\,f(\fp)$ is the center of $v$,
{\em i.e.} the image of $\fm_V$ in $\Spec\,A$
(see remark \ref{rem_Spv-of-ring}(v)). Then we may find
$\fp'\in\Spec\,B'$ whose images in $\Spec\,B$ and $\Spec\,V$
equal respectively $\fp$ and $\fm_V$. We pick a minimal
prime ideal $\fq'$ of $B'$ contained in $\fp'$. Since
$f$ is flat, the same holds for the induced map
$V\to B\otimes_AV$, and therefore the image of $\fq'$
in $\Spec\,V$ equals $\eta_V$ (\cite[Th.9.5]{Mat});
consequently the image of $\fq'$ in $\Spec\,B$ lies in
$U_B$. Let $\bar b_0\in B'$ be the image of $b_0$; from
condition (D) and \cite[Th.6.1(ii), Th.6.5(iii)]{Mat}, we
deduce that
\set\begin{equation}\label{eq_b_0-not-here}
\bar b_0\notin\fq'.
\end{equation}
Choose a valuation ring $V'$ of the residue field
$\kappa(\fq')$ that dominates the image of the local
ring $B'_{\fp'}$ (corollary \ref{cor_cornerstone}).
Then, $V'$ corresponds to a point $v'\in\Spv\,B$,
and we notice that $\Spv\,f(v')=v$ : indeed, the induced
map $V\to\kappa(\fq')$ is injective, and $V'$ dominates
the image of this map, whence the assertion. To conclude
the proof of the claim, it suffices to check that
$v'\in T\cap\Spv\,U_B$.
However, \eqref{eq_b_0-not-here} already implies that
$v'(b_0)\neq 0$. Moreover, since $V'$ contains the image
of $B$, we have $v(b_i)\leq 1$ for every $i=1,\dots,m$,
and since the image of $\fp'$ lies in the maximal ideal
of $V$, we have $v(b_i)<1$ for every $i=1,\dots,n$.
Thus, $v\in T$. Lastly, since the image of $\fq'$ in
$\Spec\,B$ lies in $U_B$, we get $v\in\Spv\,U_B$ as well,
as needed.
\end{pfclaim}
\end{proof}

\begin{remark} Our proof of theorem \ref{th_Chevalley}
employs crucially Gruson and Raynaud's technique of
flattening stratifications (\cite{Gr-Ra}). On the other
hand, Huber in \cite{Hu1} recovers this theorem rather
easily via ``elimination of quantifiers'' for the first
order theory of algebraically closed valued fields with
non-trivial valuation. However, the latter is of course
a deep model theoretic result, so Huber's proof is hardly
more elementary than ours.
\end{remark}

We wish next to explain how the valuation spectrum can
be applied to the study of integral closures of rings
and ideals. We begin with the following :

\begin{definition}\label{def_int-cl-ideal}
Let $f:A\to B$ be any ring homomorphism, $I\subset A$
any ideal. The {\em integral closure of\/ $I$ in $B$}
is the set denoted
$$
\mathrm{i.c.}(I,B)
$$
and consisting of all $x\in B$ that satisfy, for some
$n\in\N$, an equation of the type :
\set\begin{equation}\label{eq_int-equation}
x^n+f(a_1)\cdot x^{n-1}+\cdots+f(a_n)=0
\qquad
\text{where $a_j\in I^j$ for $j=1,\dots,n$}.
\end{equation}
\end{definition}

Notice that for $I=A$, the set $\mathrm{i.c.}(A,B)$
is the usual integral closure of $A$ in $B$. We have :

\begin{lemma}\label{lem_int-closure-Rees}
In the situation of definition {\em\ref{def_int-cl-ideal}},
let $b\in B$ be any element.
\begin{enumerate}
\item
The following conditions are equivalent :
\begin{enumerate}
\item
$b\in\mathrm{i.c.}(I,B)$.
\item
$bU^{-1}\in\mathrm{i.c.}(\sR(A,I)_\bullet,\sR(B,B)_\bullet)$
(notation of example {\em\ref{ex_Rees}}).
\end{enumerate}
\item
Especially, $\mathrm{i.c.}(I,B)$ is an ideal of\/
$\mathrm{i.c.}(A,B)$.
\item
For every prime ideal $\fp\in\Spec\,B$, let
$\bar b_\fp\in\kappa(\fp)$ denote the image of\/ $b$
(notation of definition {\em\ref{def_strict-loc}(i)}).
The following conditions are equivalent :
\begin{enumerate}
\item
$b\in\mathrm{i.c.}(I,B)$.
\item
$\bar b_\fp\in\mathrm{i.c.}(I,\kappa(\fp))$\ \ for every
minimal prime ideal $\fp$ of $B$.
\end{enumerate}
\end{enumerate}
\end{lemma}
\begin{proof} Suppose that $b$ satisfies an identity
\eqref{eq_int-equation}, and set $T:=U^{-1}$ to ease
notation; we get
$$
(bT)^n+Tf(a_1)\cdot(bT)^{n-1}+\cdots+T^nf(a_n)=0
$$
where $T^jf(a_j)\in\sR(A,I)_\bullet$ for $j=1,\dots,n$.
Conversely, suppose the $bT$ is integral over
$\sR(A,I)_\bullet$, so that we have a monic polynomial
$P(X)=X^n+u_1X^{n-1}+\cdots+u_n$, with coefficients in
$\sR(A,I)_\bullet$, such that $P(bT)=0$; for every
$j=1,\dots,n$, denote by $v_j$ the component of $u_j$
in degree $n$ (for the grading of $\sR(A,I)_\bullet$);
since $bT$ is homogeneous of degre $1$, it follows
that $(bT)^n+v_1\cdot(bT)^{n-1}+\cdots+v_n=0$, and
letting $T=1$, we deduce that $b$ is integral over $I$.
This proves (i), and assertion (ii) is an immediate
consequence.

(iii): Clearly (iii.a)$\Rightarrow$(iii.b). Hence,
suppose that (iii.b) holds, and for every minimal
prime ideal $\fp$ of $B$, pick a monic polynomial
$P_\fp(X)\in A[X]$, say of degree $n_\fp$, such that
$P_\fp(\bar b_\fp)=0$, and whose coefficient in
degree $n_\fp-j$ lies in $I^j$, for $j=1,\dots n_\fp$.
For every such $\fp$, let $V_\fp\subset\Spec\,B$ denote
the subset of all prime ideals containing $P_\fp(b)$.
Hence, $Z:=\bigcup_{\fp\in\Min(B)}V_\fp$ contains the set
$\Min(B)$ of minimal prime ideals of $B$; but each
$V_\fp$ is a constructible closed subset of $\Spec\,B$,
hence $Z=\Spec\,B$, and there exists a finite subset
$S\subset\Min(B)$ such that
$\bigcup_{\fp\in S}V_\fp=\Spec\,B$ (theorem
\ref{th_main-spectral}(i)). Set $P:=\prod_{\fp\in S}P_\fp$;
it follows that $P(b)$ lies in the nilpotent ideal of
$B$, so $P^r(b)=0$ for a sufficiently large $r\in\N$.
However, let $N$ be the degree of $P^r$; it is easily
seen that the coefficient in degree $N-j$ of $P^r$ lies
in $I^j$, for every $j=1,\dots,N$ (details left to the
reader), whence (iii.a).
\end{proof}

Concerning integral closures of ideals, much more can
be found in the treatise \cite{Hu-Sw}, a comprehensive
reference for this subject. For our purposes, proposition
\ref{prop_integ-clos-I} hereafter will suffice. For its
proof, we shall employ the technical lemma
\ref{lem_technical-lemma}, which shall be applied also
later to the study of the adic spectrum of a topological
ring, and for this reason, is stated here in wider
generality than needed for our immediate purposes in
this section.

\sset\subsubsection{}\label{subsec_technical}
Let $B$ be any ring, $A\subset B$ a subring, $I$ an ideal
of $A$, and $\Sigma,\Sigma'$ two sets; denote by $\cP_I$
the subset of the free polynomial $A$-algebra
$$
R[Z]
\qquad\text{where}\qquad
R:=A[X_i,Y_j~|~i\in\Sigma,\ j\in\Sigma']
$$
consisting of all polynomials $P$ fulfilling the
following conditions :
\begin{itemize}
\item
We have $P=Z^n+\sum_{k=1}^nZ^{n-k}\cdot P_k(X_\bullet,Y_\bullet)$
for some $P_1,\dots,P_n\in R$ such that $P_k$ is homogeneous
of degree $k$ for every $k=1,\dots,n$ (for the
standard $\N$-grading on $R$ such that $X_i,Y_j\in\deg_1R$
for every $i\in\Sigma$ and $j\in\Sigma'$).
\item
For every $k=1,\dots,n$, denote by $P_k(0,Y_\bullet)$
the image of $P_k$ by the evaluation map
$R\to A[Y_j~|~j\in\Sigma']$ such that $X_i\mapsto 0$
for every $i\in\Sigma$; then
$P_k(0,Y_\bullet)\in I[Y_i~|~j\in\Sigma']$.
\end{itemize}
Set as well
$$
\begin{aligned}
T:= &\, \{v\in\Spv\,B~|~\text{$v(a)\leq 1$ and $v(b)<1$
for every $a\in A$ and $b\in I$}\} \\
T_0:= &\, 
\{v\in T~|~\text{$\Ker\,v$ is a minimal prime ideal of $B$}\}.
\end{aligned}
$$
Notice that $T$ is a pro-constructible subset of $\Spv\,B$,
and it contains all the primary generizations of its points.

\begin{lemma}\label{lem_technical-lemma}
In the situation of \eqref{subsec_technical}, suppose
moreover that $I$ is contained in the Jacobson radical
of $A$, and let also $(g_i~|~i\in\Sigma)$ and
$(h_j~|~j\in\Sigma')$ be two systems of elements of $B$.
Then for every $f\in B$ the following conditions are
equivalent :
\begin{enumerate}
\alphaenu
\item
For every $v\in T$ with $v(f)\neq 0$ there exists either
$i\in\Sigma$ such that $v(f)\leq v(g_i)$ or else $j\in\Sigma'$
such that $v(f)<v(h_j)$.
\item
For every $v\in T_0$ with $v(f)\neq 0$ there exists either
$i\in\Sigma$ such that $v(f)\leq v(g_i)$ or else $j\in\Sigma'$
such that $v(f)<v(h_j)$.
\item
There exists a polynomial $P\in\cP_I$ such that
$P(f,g_\bullet,h_\bullet)=0$ in $B$.
\end{enumerate}
\end{lemma}
\begin{proof}(c)$\Rightarrow$(a): Let $P\in\cP_I$ such that
(c) holds, and write
$P=Z^n+\sum_{k=1}^nZ^{n-k}\cdot P_k(X_\bullet,Y_\bullet)$ for
some $P_1,\dots,P_n\in R$. Fix $v\in T$; then there exists
$k\leq n$ such that
$$
v(f^n)\leq v(f^{n-k}\cdot P_k(g_\bullet,h_\bullet)).
$$
If $v(f)=0$, there is nothing to show; otherwise, we deduce
that $v(f^k)\leq v(P_k(g_\bullet,h_\bullet))$. By assumption,
$P_k$ is homogeneous of degree $k$, so there exist two
systems of non-negative integers
$\nu_\bullet:=(\nu_i~|~i\in\Sigma)$,
$\mu_\bullet:=(\mu_j~|~j\in\Sigma')$, and an element
$a\in A$ such that the monomial
$aX_\bullet^{\nu_\bullet}Y_\bullet^{\mu_\bullet}$ appears in $P_k$
and
$$
v(f^k)\leq
v(a)\cdot v(g_\bullet^{\nu_\bullet})\cdot v(h_\bullet^{\mu_\bullet}).
$$
Especially, set $|\nu_\bullet|:=\sum_{i\in\Sigma}\nu_i$, and
define likewise $|\mu_\bullet|$; then
$|\nu_\bullet|+|\mu_\bullet|=k$. Suppose first that
$|\nu_\bullet|>0$; if there exists $i\in I$ with
$v(f)\leq v(g_i)$, we are done. Otherwise, we have
$v(f^{|\nu_\bullet|})>v(g_\bullet^{\nu_\bullet})$, and therefore
$v(f^{|\mu_\bullet|})<v(h_\bullet^{\mu_\bullet})$, since $v(a)\leq 1$.
In this case, we easily get $v(f)<v(h_j)$ for
some $j\in\Sigma'$. Lastly, if $|\nu_\bullet|=0$, by
assumption we must have $a\in I$, so that $v(a)<1$,
and therefore $v(f^k)<v(h_\bullet^{\mu_\bullet})$, so we
conclude as in the previous case.

(a)$\Rightarrow$(b) is obvious.

(b)$\Rightarrow$(c): We suppose that (c) fails, and we
exhibit $v\in T_0$ such that $v(f)\neq 0$ and $v(g_i)<v(f)$,
$v(h_j)\leq v(f)$ for every $i\in\Sigma$ and $j\in\Sigma'$.
To this aim, consider the subring
$$
A':=A[g_i/f,h_j/f~|~i\in\Sigma,\ j\in\Sigma']
\subset B':=B[f^{-1}]
$$
and the ideal $J_0$ of $A'$ generated by the
system $(g_i/f~|~i\in\Sigma)$; set $J:=J_0+IA'$.

\begin{claim}\label{cl_if-b-fails}
If (c) fails, we have $J\neq A'$.
\end{claim}
\begin{pfclaim} If $J=A'$, we may find an expression
of the form
$$
1=\sum_{i\in\Sigma}\frac{g_i}{f}\cdot
F_i\Bigl(\frac{g_\bullet}{f},\frac{h_\bullet}{f}\Bigr)+
\sum_{k=1}^ra_k\cdot
G_k\Bigl(\frac{g_\bullet}{f},\frac{h_\bullet}{f}\Bigr)
\qquad
\text{in $A'$}
$$
for systems $(F_i(X_\bullet,Y_\bullet)~|~i\in\Sigma)$ and
$(G_k(X_\bullet,Y_\bullet)~|~k=1,\dots,r)$ of elements of $R$,
with $F_i=0$ except for finitely many values of $i\in\Sigma$,
and a system $(a_1,\dots,a_r)$ of elements of $I$.
We may then multiply both sides of this identity by $f^{n+1}$
for a sufficiently large $n\in\N$, and obtain an identity
in $A$ of the form $P(f,g_\bullet,h_\bullet)=0$, with
$$
P(Z,X_\bullet,Y_\bullet)=Z^{n+1}-
\sum_{i\in\Sigma}Z^n\cdot X_i\cdot F_i(X_\bullet/Z,Y_\bullet/Z)
-\sum_{k=1}^r Z^{n+1}\cdot a_k\cdot G_k(X_\bullet/Z,Y_\bullet/Z).
$$
A simple inspection shows that the coefficient of the
monomial $Z^{n+1}$ appearing in $P$ is of the form
$1-a$ for some $a\in I$ : namely, $a$ is the sum
of the terms $a_k$ such that $G_k\in A$. However,
$1-a\in A^\times$, since $I$ is contained in the Jacobson
radical of $A$; to conclude the proof, it suffices to
observe that $(1-a)^{-1}\cdot P\in\cP_I$ : details left to
the reader.
\end{pfclaim}

By virtue of claim \ref{cl_if-b-fails}, we may find
a maximal ideal $\fm$ of $A'$ containing $J$, and we
let $\fp\subset\fm$ be any minimal prime ideal of $A'$;
then $A'_\fp\subset B'_\fp$ and we pick any minimal prime
ideal $\fq$ of $B'_\fp$. Clearly $\fq\cap A'_\fp=\fp A'_\fp$,
so the induced map $C:=A'/\fp\to D:=B'_\fp/\fq$ is
still injective. Let $\bar\fm\subset C$ be the image
of $\fm$; by corollary \ref{cor_cornerstone}, we may
find a valuation ring $V$ of $\Frac\,D$ that dominates
$C_{\bar\fm}$. The valuation ring $V$ corresponds to an
element $v\in T_0$ with the required properties.
\end{proof}

\begin{proposition}\label{prop_integ-clos-I}
Let $i:A\to B$ be an injective ring homomorphism,
$J$ any ideal of $A$. Set
$$
\cV_{B/A}:=\Spv(i)^{-1}(\Spv^+A)
\qquad
\cV^0_{B/A}:=\{v\in\cV_{B/A}~|~
\text{$\ker\,v$ is a minimal prime ideal of $B$}\}.
$$
For every $v\in\cV_{B/A}$, let also $\pi_v:B\to\kappa(v)$
be the natural map. Then we have :
$$
\mathrm{i.c.}(J,B)=\!\!\!
\bigcap_{v\in\cV_{B/A}}\!\!\!\pi_v^{-1}(J\cdot\kappa(v)^+)=
\!\!\!\bigcap_{v\in\cV^0_{B/A}}\!\!\!\pi_v^{-1}(J\cdot\kappa(v)^+).
$$
\end{proposition}
\begin{proof} By applying lemma \ref{lem_technical-lemma}
with $I=0$, $(g_i~|~i\in\Sigma)$ any (non-empty) set of
generators for $J$, and $\Sigma'=\emptyset$, we see that
$\mathrm{i.c.}(J,B)$ is the set of elements $f\in B$
such that the following holds. For every $v\in\cV_{B/A}$
(respectively, for every $v\in\cV^0_{B/A}$) there exists
$i\in\Sigma$ with $v(f)\leq v(g_i)$. The proposition is
an immediate consequence.
\end{proof}

\begin{remark}\label{rem_add-a-parameter}
For every ring homomorphism $f:A\to B$, every ideal
$I\subset A$, and every $q\in\Q_+$, write $q=m/n$ for
integers $m,n\in\N$, and let
$$
\mathrm{i.c.}(I,B,q):=\{b\in B~|~b^n\in\mathrm{i.c.}(I^m,B)\}.
$$

(i)\ \
Define $\cV_{B/A}:=\Spv(f)^{-1}(\Spv^+A)$.
We claim that $\mathrm{i.c.}(I,B,q)$ is the set of all
$b\in B$ fulfilling the following condition. For every
$v\in\cV_{B/A}$ there exists $x\in I$ such that
$v(b)\leq v(f(x))^q$ (see remark \ref{rem_ordered-gps}(ii)).
Indeed, by proposition \ref{prop_integ-clos-I}, we have :
$$
\mathrm{i.c.}(I,B,q)=
\{b\in B~|~\text{$\pi_v(b^n)\in I^m\cdot\kappa(v)^+$
for every $v\in\cV_{B/A}$}\}.
$$
The latter condition can be restated as follows. For every
$v\in\cV_{B/A}$ there exists $x\in I^m\cdot\kappa(v)^+$
such that $v(b^n)\leq v(x)$. We may write such $x$ as a
sum $\sum_{i=1}^kf(x_i)\cdot y_i$ for certain $x_i\in I^m$,
$y_i\in V$, and then $v(x)\leq v(f(x_i))$ for at least an
index $i\leq r$. By the same token, we may assume that
$x_i$ is a product $a_1a_2\cdots a_m$ of elements of $I$,
and then $v(x_i)\leq v(f(a_j)^m)$ for at least one index
$j\leq m$.

(ii)\ \
Especially, (i) shows that $\mathrm{i.c.}(I,B,q)$ depends
only on $q$ (and not on the integers $m,n$). 
Notice as well that
\set\begin{equation}\label{eq_lower-m}
\mathrm{i.c.}(I^m,B)=\mathrm{i.c.}(I,B,m)
\qquad
\text{for every $m\in\N$}.
\end{equation}
Moreover, a simple inspection of the definition reveals that
$\mathrm{i.c.}(I,B,q)\subset\mathrm{i.c.}(A,B)$, and then
the criterion of (i) also implies that $\mathrm{i.c.}(I,B,q)$
is an ideal of $\mathrm{i.c.}(A,B)$, for every $q\in\Q_+$.
Furthermore, we have
\set\begin{equation}\label{eq_twice-ic}
\mathrm{i.c.}(I,B,q)=\mathrm{i.c.}(\mathrm{i.c.}(I,A),B,q)
\qquad
\text{for every $q\in\Q_+$}.
\end{equation}
Indeed, since $I\subset\mathrm{i.c.}(I,A)$, we have obviously
$\mathrm{i.c.}(I,B,q)\subset\mathrm{i.c.}(\mathrm{i.c.}(I,A),B,q)$.
Conversely, if $b\in\mathrm{i.c.}(\mathrm{i.c.}(I,A),B,q)$ and
$v\in\cV_{B/A}$ is any element, (i) says that there exists
$x\in\mathrm{i.c.}(I,A)$ such that $v(b)\leq v(f(x))^q$;
by the same token, there exists $y\in I$ such that
$v(x)\leq v(y)$, so that $v(b)\leq v(f(y))^q$, which shows
that $b\in\mathrm{i.c.}(I,B,q)$. Next, let us set
$$
\bar{\mathrm{i.c.}}(I,B,r):=
\bigcap_{\substack{ q\in\Q_+ \\
                  q<r\ \ }}
\mathrm{i.c.}(I,B,q)
\qquad
\text{for every $r\in\R_{>0}$}
$$
and define as well
$\bar{\mathrm{i.c.}}(I,B,0):=\mathrm{i.c.}(A,B)$. 
Denote also by $\cR_{B/A}$ the set of all real valued
valuations $|\cdot|$ on $B$ that lie in $\cV_{B/A}$.
We have :
\end{remark}

\begin{lemma}\label{lem_strict-closure}
With the notation of remark
{\em\ref{rem_add-a-parameter}(ii)}, let $r\in\R_+$ be
any real number, $b\in B$ any element, and denote
by $A[b]$ the $A$-subalgebra of $B$ generated by $b$.
The following conditions are equivalent :
\begin{enumerate}
\alphaenu
\item
$b\in\bar{\mathrm{i.c.}}(I,B,r)$.
\item
$|b|\leq\sup_{x\in I}|f(x)|^r$\ \ for every element $|\cdot|$
of\/ $\cR_{A[b]/A}$.
\end{enumerate}
\end{lemma}
\begin{proof}(a)$\Rightarrow$(b): The assumption means
that $b\in\mathrm{i.c.}(I,A[b],q)$ for every non-negative
rational number $q<r$. By remark \ref{rem_add-a-parameter}(i),
this implies that, for every $|\cdot|$ in $\cR_{A[b]/A}$
and every such $q$, there exists $x\in I$ such that
$|b|\leq|f(x)|^q$. Assertion (b) is an immediate consequence.

(b)$\Rightarrow$(a): Suppose first that $r=0$. In this
case, we have $|b|\leq 1$ for every $\cR_{A[b]/A}$, and
we need to show that $b\in\mathrm{i.c.}(A,B)$.
Suppose, that the latter fails; then, by remark
\ref{rem_add-a-parameter}(i), there exists a valuation
$|\cdot|$ in $\cV_{A[b]/A}$ such that $|b|>1$.
It is easily seen that the characteristic subgroup of
$|\cdot|$ equals $\Delta:=O(|b|)$ (notation
of \eqref{subsec_special-investig} and remark
\ref{rem_ordered-gps}(vii)), hence $|\cdot|^\Delta$ is
still an element of $\cV_{A[b]/A}$ and $|b|^\Delta>1$
(lemma \ref{lem_specialize-vals}(ii)). We may therefore
assume that the value group of $|\cdot|$ equals $\Delta$,
and $|\cdot|':=|\cdot|_{o(|b|)}$ is then a valuation of
rank one (lemma \ref{lem_specialize-vals}(i) and remark
\ref{rem_ordered-gps}(vii)). Lastly, $|\cdot|'$
is equivalent to an element of $\cR_{A[b]/A}$, by remark
\ref{rem_ordered-gps}(vi), and by construction we have
$|b|'>1$, a contradiction.

Next, suppose that $q>0$, in which case -- by remark
\ref{rem_add-a-parameter}(i) -- it suffices to show
that, for every $|\cdot|$ in $\cV_{A[b]/A}$ and every
non-negative rational number $q<r$, there exists $x\in I$
such that $|b|\leq|f(x)|^q$. However, by the foregoing case
we know already that $|b|\leq 1$ for every such $|\cdot|$. 

$\bullet$\ \ 
If $|b|=0$, there is nothing to prove, so suppose now
that $0<|b|<1$, and set $\Delta:=O(|b|)$, $\Delta':=o(|b|)$;
notice also that the characteristic subgroup of $|\cdot|$
is $\{1\}$, hence $|\cdot|^\Delta$ lies still in $\cV_{A[b]/A}$
(lemma \ref{lem_specialize-vals}(ii)) and $|b|^\Delta=|b|$.
The same then holds for $|\cdot|':=|\cdot|^\Delta_{\Delta'}$
(lemma \ref{lem_specialize-vals}(i)), and as in the
foregoing, we see that the latter valuation is equivalent
to an element of $\cR_{A[b]/A}$. Moreover, since
$|b|\notin\Delta'$, we still have $|b|'<1$. Assumption (b)
then implies that $|b|'\leq\sup_{x\in I}|f(x)|'^r$. Explicitly,
for every real number $\eps>0$ there exists $x_\eps\in I$
such that $|f(x_\eps)|'^r>|b|'-\eps$. Then
$|f(x_\eps)|'^q>t_\eps:=(|b|'-\eps)^{r/q}$, and since $q/r<1$
and $|b|'<1$, we may find $\eps$ small enough, so that
$t_\eps>|b|'$, and consequently $|f(x_\eps)|'^q>|b|'$.
It follows that $|f(x_\eps)|^q>|b|$, as required.

$\bullet$\ \
Lastly, suppose that $|b|=1$. In this case, we must exhibit
some $x\in I$ with $|f(x)|=1$. Suppose that no such element
exist, and let $\pi:A[b]\to\kappa:=\kappa(|\cdot|)$ be the
projection (notation of remark \ref{rem_semi-norm}(v)); in the
current situation, the image of $\pi$ lies in the valuation
ring $\kappa^+$ of the induced valuation $|\cdot|_\kappa$,
and the assumption means that $\pi(I)$ lies in the maximal
ideal $\fm$ of $\kappa^+$. Then, denote by $|\cdot|'_\kappa$
the trivial valuation on $\kappa^+/\fm$; the composition of
$|\cdot|'_\kappa$ with the induced projection $A\to\kappa^+/\fm$
gives an element $|\cdot|'_A$ of $\cR_{A[b]/A}$ such that
$|b|'_A=1$ and $|f(x)|'_A=0$ for every $x\in I$, contradicting
assumption (b).
\end{proof}

For future reference, we point out the following
``approximation lemma'' for integral closures:

\begin{lemma}\label{lem_approx-ic}
Let $A$ be any ring, $a_\bullet:=(a_1,\dots,a_n)$ and
$b_\bullet:=(b_1,\dots,b_n)$ two sequences of elements
of $A$, and $q>1$ any rational number such that
$$
b_i-a_i\in\mathrm{i.c.}(I,A,q)
\qquad
\text{for every $i=1,\dots,n$}.
$$
Denote by $I\subset A$ (resp. $J\subset A$) the ideal
generated by $a_\bullet$ (resp. by $b_\bullet$), and suppose
that :
\begin{enumerate}
\alphaenu
\item
either, the radical of $I$ equals the radical of $J$
\item
or else, $I$ is contained in the Jacobson radical of $A$.
\end{enumerate}
Then $\mathrm{i.c.}(I,A)=\mathrm{i.c.}(J,A)$.
\end{lemma}
\begin{proof}Since $\mathrm{i.c.}(I,A,q)\subset\mathrm{i.c.}(I,A)$,
it is clear that $\mathrm{i.c.}(J,A)\subset\mathrm{i.c.}(I,A)$.
To show the converse inclusion, set
$$
r_v:=\max(v(a_1),\dots,v(a_n))
\qquad
s_v:=\max(v(b_1),\dots,v(b_n))
\qquad
\text{for every $v\in\Spv\,A^+$}.
$$
By virtue of proposition \ref{prop_integ-clos-I}, it suffices
to check that $s_v\geq r_v$ for every such $v$. Now, suppose
first that (a) holds, and say that $r_v=v(a_i)$ for some
$i\leq n$; under our assumptions, remark
\ref{rem_add-a-parameter}(i) implies that
$v(b_j-a_j)\leq r_v^q$ for every $j=1,\dots,n$. If $r_v=1$,
the ideal $I$ is not contained in the center of $v$ (see
remark \ref{rem_Spv-of-ring}(v)), and then the same holds
for $J$, so that $s_v=1$ as well. If $r_v=0$, we get
$v(b_j-a_j)=0$ for every $j\leq n$, whence $s_v=0$ as
well. Lastly, if $0<s_v<1$, we get $s_v^q<s_v$, so that
$v(a_i)=v(b_i)$, and thus $s_v\geq r_v$ again.

Next, suppose that (b) holds, and let $v\in\Spv^+A$ be any
valuation; again, we need to show that $s_v\geq r_v$, and by
virtue of proposition \ref{prop_prprties-of-sigma}(iii) we
may assume that the center $\fp$ of $v$ is a maximal ideal
of $A$; especially $I\subset\fp$, in which case $r_v<1$,
and we argue as in the previous case, to conclude.
\end{proof}

We shall continue in section \ref{sec_Witt-Fontaine}
the study of the class of ideals considered in remark
\ref{rem_add-a-parameter}. We shall be especially
interested in the case where $A=B$ is a perfect
topological $\F_p$-algebra: see lemma
\ref{lem_return-to-ic}, which in turn shall be useful for
our investigation of the finer topological properties of
the ring $W(A)$ of Witt vectors over $A$.

\subsection{Witt vectors}
\label{sec_Witt-Fontaine}
We begin with a review of the ring of $p$-typical Witt
vectors associated with an arbitrary ring; the basic
reference is \cite[Ch.IX, \S1]{BouAC}, though our treatment
is self-contained and pays greater attention to the
topological aspects of the theory. The second topic
of this section is the study of the functor here denoted
$\bE$ (see \eqref{subsec_ring-struct-on-E}), that
generalizes Fontaine and Wintenberger's {\em field of norms}
of a local field of characteristic zero. In section
\ref{sec_def-A-tilde}, these two constructions will
merge in the definition of the functor $\bA$, that will
play an important role in the study of perfectoid rings.

Henceforth we fix a prime number $p$ and we set
$\F_p:=\Z/p\Z$.

\begin{definition}
Let $M$ be the free monoid with basis $\{X_i,Y_i~|~i\in\N\}$,
and consider the morphism of monoids
$$
\deg:M\to\N^{\oplus 2}
\qquad\text{such that}\quad
\deg(X_i):=(p^i,0)
\quad
\deg(Y_i):=(0,p^i)
\qquad
\text{for every $i\in\N$}.
$$
Let $P\in\Z[M]=\Z[X_i,Y_i~|~i\in\N]$ be any polynomial;
we can write uniquely $P=\sum^k_{i=1}n_iP_i$ for integers
$n_1,\dots,n_k\in\Z\setminus\{0\}$ and distinct monomials
$P_1,\dots,P_k\in M$, and we say that :
\begin{itemize}
\item
$P$ is {\em bihomogeneous of degree $(a,b)\in\N^{\oplus 2}$},
if $\deg(P_i)=(a,b)$ for every $i=1,\dots,k$
\item
$P$ is {\em homogeneous of total degree $d\in\N$}, if
$\deg(P_i)=(a_i,b_i)$ with integers $a_i,b_i\in\N$ such
that $a_i+b_i=d$, for every $i=1,\dots,k$.
\end{itemize}
\end{definition}

\sset\subsubsection{}\label{subsec_Witt-ghost}
Following \cite[Ch.IX, \S1, n.1]{BouAC} one defines,
for every $n\in\N$, the {\em $n$-th Witt polynomial}
$$
\bomega_n(X_0,\dots,X_n):=
\sum^n_{i=0}p^iX_i^{p^{n-i}}\in\Z[X_0,\dots,X_n].
$$
Notice the inductive relations :
\set\begin{equation}\label{eq_witt-pols}
\begin{aligned}
\bomega_{n+1}(X_0,\dots,X_{n+1}) =\: &
\bomega_n(X_0^p,\dots,X^p_n)+p^{n+1}X_{n+1} \\
=\: & X_0^{p^{n+1}}+p\cdot\bomega_n(X_1,\dots,X_{n+1}).
\end{aligned}
\end{equation}
Let $A$ be any ring; for every $n\in\N$, and every
$\underline a:=(a_n~|~n\in\N)\in A^\N$, the element
$$
\bomega_n(\underline a):=\bomega_n(a_0,\dots,a_n)\in A
$$
is called the {\em $n$-th ghost component\/} of
$\underline a$, and we consider the {\em ghost map}
$$
\bomega_A:A^\N\to A^\N
\qquad
\underline a\mapsto(\bomega_n(\underline a)~|~n\in\N).
$$

\begin{lemma}\label{lem_basic-cong}
Let $A$ be a ring, $(J_n~|~n\in\N)$ a sequence of ideals
of $A$, such that
$$
pJ_n+J^p_n\subset J_{n+1}\subset J_n
\qquad
\text{for every $n\in\N$}.
$$
Let also $x,y\in A$ be any two elements, and
$\underline a,\underline b\in A^\N$ any two sequences.
The following holds :
\begin{enumerate}
\item
If $x\equiv y \pmod{J_0}$, we have
$x^{p^n}\equiv y^{p^n} \pmod{J_n}$ for every $n\in\N$.
\item
If $a_i\equiv b_i \pmod{J_m}$ for every $i=0,\dots,n$,
then $\bomega_i(\underline a)\equiv\bomega_i(\underline b)
\pmod{J_{m+i}}$ for every $i=0,\dots,n$.
\end{enumerate}
\end{lemma}
\begin{proof}(i) is left to the reader, and (ii) follows
easily from (i).
\end{proof}

\begin{proposition}\label{prop_basic-cong}
With the notation of \eqref{subsec_Witt-ghost}, we have :
\begin{enumerate}
\item
If $p$ is a regular (resp. invertible) element in $A$,
then $\bomega_A$ is injective (resp. bijective).
\item
Suppose that $A$ admits a ring endomorphism $\phi:A\to A$
such that
$$
\phi(a)\equiv  a^p \pmod{pA}
\qquad
\text{for every $a\in A$}.
$$
Then $\Img\,\bomega_A=\{\underline b\in A^\N~|~
\phi(b_n)\equiv b_{n+1} \pmod{p^{n+1}A}\ \ 
\text{for every $n\in\N$}\}$.
Especially, $\Img\,\bomega_A$ is a subring of $A^\N$.
\end{enumerate}
\end{proposition}
\begin{proof}(i): Indeed, \eqref{eq_witt-pols} shows that
the condition $\bomega_n(\underline a)=\underline b$ is
equivalent to the identities
$$
b_0=a_0
\qquad\text{and}\qquad
b_n=\bomega_{n-1}(a_0^p,\dots,a_{n-1}^p)+p^na_n
\qquad
\text{for every $n\geq 1$}
$$
whence the assertion.

(ii): An easy induction argument reduces the assertion
to the following :

\begin{claim} Let $n>0$ be any integer, and
$a_0,\dots,a_{n-1},b_n\in A$ any sequence of elements.
Set $b_{n-1}:=\bomega_{n-1}(a_0,\dots,a_{n-1})$. Then
the following conditions are equivalent :
\begin{enumerate}
\alphaenu
\item
There exists $a_n\in A$ such that
$b_n=\bomega_n(a_0,\dots,a_n)$.
\item
$\phi(b_{n-1})\equiv b_n \pmod{p^nA}$.
\end{enumerate}
\end{claim}
\begin{pfclaim}[] We apply lemma \ref{lem_basic-cong}(ii)
to the filtration $(p^nA~|~n\in\N)$, to deduce that
$$
\begin{aligned}
\bomega_n(a_0,\dots,a_{n-1},x)\equiv\, &
\bomega_{n-1}(a_0^p,\dots,a^p_{n-1}) & \pmod{p^nA} &
\qquad\text{(by \eqref{eq_witt-pols})} \\
\equiv\, & \bomega_{n-1}(\phi(a_0),\dots,\phi(a_{n-1}))
& \pmod{p^nA} \\
\equiv\, & \phi(b_{n-1}) & \pmod{p^nA}
\end{aligned}
$$
for every $x\in A$. The claim is an immediate consequence.
\end{pfclaim}
\end{proof}

\sset\subsubsection{}\label{subsec_Witt-laws}
Take now $A:=\Z[X_i,Y_i~|~i\in\N]$, and let $\phi:A\to A$
be the endomorphism such that $\phi(X_i)=X_i^p$ and
$\phi(Y_i)=Y_i^p$ for every $i\in\N$. Clearly $\phi$
fulfills the condition of proposition \ref{prop_basic-cong}(ii);
from proposition \ref{prop_basic-cong} and
\eqref{eq_witt-pols} we deduce that there exist polynomials
$$
S_n,P_n\in\Z[X_0,\dots,X_n,Y_0,\dots,Y_n]
\quad
I_n\in\Z[X_0,\dots,X_n]
\quad
F_n\in\Z[X_0,\dots,X_{n+1}]
$$
uniquely characterized, for every $n\in\N$, by the identities :
$$
\begin{aligned}
\bomega_n(S_0,\dots,S_n) =\: &
\bomega_n(X_0,\dots,X_n)+\bomega_n(Y_0,\dots,Y_n) \\
\bomega_n(P_0,\dots,P_n) =\: &
\bomega_n(X_0,\dots,X_n)\cdot\bomega_n(Y_0,\dots,Y_n) \\
\bomega_n(I_0,\dots,I_n) =\: & -\bomega_n(X_0,\dots,X_n) \\
\bomega_n(F_0,\dots,F_n) =\: & \bomega_{n+1}(X_0,\dots,X_{n+1}).
\end{aligned}
$$

\begin{remark}\label{rem_homogeneous-laws}
(i)\ \
Notice that for every $n\in\N$, the polynomial
$\omega_n$ is bihomogeneous of degree $(p^n,0)$.
By a simple induction argument, it follows easily
that $S_n$ is homogeneous of total degree $p^n$,
and $P_n$ (resp. $I_n$, resp. $F_n$) is bihomogeneous
(resp. homogeneous) of degree $(p^n,p^n)$ (resp. $p^n$,
resp. $p^{n+1}$), for every $n\in\N$.

(ii)\ \
For example, a simple calculation yields the formulae :
$$
\begin{aligned}
S_0=\: & X_0+Y_0 \qquad & S_1=\: &
X_1+Y_1-\sum_{i=1}^{p-1}\frac{1}{p}\binom{p}{i}X_0^iY_0^{p-i} \\
P_0=\: & X_0Y_0 \qquad & P_1=\: & X_0^pY_1+X_1Y_0^p+pX_1Y_1 \\
I_0=\: & -X_0 \qquad & I_1=\: & -X_1-\eps\cdot X_0^p \\
F_0=\: & X_0^p+pX_1 \qquad &
F_1=\: & (1-p^{p-1})X_1^p+pX_2-
         \sum_{i=1}^{p-1}\frac{1}{p}\binom{p}{i}X_0^{pi}(pX_1)^{p-i}. \\
\end{aligned}
$$
where $\eps=1$ if $p=2$, and $\eps=0$ if $p>2$. Also, since
$\omega_n(0,\dots,0,X_n)=p^nX$, we see that
\set\begin{equation}\label{eq_addition-on-Vn}
S_n(0,\dots,0,X_n,0,\dots,0,Y_n)=X_n+Y_n
\qquad
\text{for every $n\in\N$}.
\end{equation}

(iii)\ \
A simple induction shows that :
\set\begin{equation}\label{eq_congr-for-F_n}
F_n\equiv X^p_n \pmod{p\Z[X_0,\dots,X_{n+1}]}
\qquad\text{for every $n\in\N$.}
\end{equation}
Indeed, for $n=0$ the assertion is clear from (ii).
Suppose that $n>0$, and that the assertion is already known
for every integer $<n$, and to ease notation, set
$A:=\Z[X_0,\dots,X_{n+1}]$ and
$\bomega_n(F):=\bomega_n(F_0,\dots,F_n)$; from
\eqref{eq_witt-pols} we get
$$
\bomega_n(F)\equiv
\bomega_n(X^p_0,\dots,X^p_n)\equiv
\bomega_{n-1}(X^{p^2}_0,\dots,X^{p^2}_{n-1})+p^nX_n^p
\pmod{p^{n+1}A}
$$
On the other hand, the inductive assumption, together
with lemma \ref{lem_basic-cong}(i,ii) implies that
$$
\bomega_n(F)\equiv
\bomega_{n-1}(F_0^p,\dots,F_{n-1}^p)+p^nF_n\equiv
\bomega_{n-1}(X_0^{p^2},\dots,X_{n-1}^{p^2})+p^nF_n
\pmod{p^{n+1}A}
$$
therefore $p^nF_n\equiv p^nX^p_n\pmod{p^{n+1}A}$,
and the assertion follows.
\end{remark}

\sset\subsubsection{}\label{subsec_Witt-vectors}
Let now $(A,+,\cdot,\cT)$ be an arbitrary topological
ring (see definition \ref{def_top-ring}(i)). We endow
$A^\N$ with the topology $\cT_{W(A)}$ defined as the product
of the topologies $\cT$ on each copy of $A$, and for every
$\underline a:=(a_n~|~n\in\N),
\underline b:=(b_n~|~n\in\N)$ in $A^\N$, we set :
$$
\begin{aligned}
S_A(\underline a,\underline b) :=\: &
(S_n(a_0,\dots,a_n,b_0,\dots,b_n)~|~n\in\N) \\
P_A(\underline a,\underline b) :=\: &
(P_n(a_0,\dots,a_n,b_0,\dots,b_n)~|~n\in\N) \\
I_A(\underline a) :=\: & (I_n(a_0,\dots,a_n)~|~n\in\N).
\end{aligned}
$$
There follow identities :
$$
\begin{aligned}
\bomega_A(S_A(\underline a,\underline b)) =\: &
\bomega_A(\underline a)+\bomega_A(\underline b) \\
\bomega_A(P_A(\underline a,\underline b)) =\: &
\bomega_A(\underline a)\cdot\bomega_A(\underline b) \\
\bomega_A(I_A(\underline a)) =\: & -\bomega_A(\underline a)
\end{aligned}
$$
for every $\underline a,\underline b\in A^\N$. Moreover,
for any other topological ring $B$, and any continuous
ring homomorphism $\psi:B\to A$, let $\psi^\N:B^\N\to A^\N$
be the induced continuous map; clearly we have the identities
\set\begin{equation}\label{eq_funct-Witt}
S_A\circ(\psi^\N\times\psi^\N)=\psi^\N\circ S_B
\quad
P_A\circ(\psi^\N\times\psi^\N)=\psi^\N\circ P_B
\quad
I_A\circ\psi^\N=\psi^\N\circ I_B.
\end{equation}

\begin{theorem}\label{th_Witt}
With the notation of \eqref{subsec_Witt-vectors}, the set
$A^\N$, endowed with the addition law
$S_A:A^\N\times A^\N\to A^\N$ and product law
$P_A:A^\N\times A^\N\to A^\N$ is a commutative topological
ring, whose zero element and unit element are respectively
the sequences
$$
\underline 0_A:=(0,0,\dots)
\qquad\text{and}\qquad
\underline 1_A:=(1,0,0,\dots).
$$
Moreover, the opposite $-\underline a$ (that is, with
respect to the addition law $S_A$) of any
$\underline a\in A^\N$, is the element $I_A(\underline a)$.
\end{theorem}
\begin{proof} It is clear that $S_A$, $P_A$ and $I_A$
are continuous mapping for the topology $\cT_{W(A)}$,
hence it remains only to prove that $(A^\N,S_A,P_A)$
is a ring. Pick any surjective ring homomorphism
$\psi:B\to A$, with $B=\Z[M]$ for some free monoid
$M$, and endow $B$ with the discrete topology; then
$\psi^\N$ is also surjective, and in light of
\eqref{eq_funct-Witt}, it suffices to prove the
theorem for $B$. Also, the $p$-Frobenius endomorphism
of $M$ induces an endomorphism $\phi$ of $B$ that
lifts the Frobenius endomorphism of $B/pB$. We may
then replace $A$ by $B$, and assume from start that
$p$ is a regular element of $A$, and $\phi:A\to A$
is a given lift of the Frobenius endomorphism of
$A/pA$. In this case, proposition
\ref{prop_basic-cong}(i,ii) implies that $\bomega_A$
maps bijectively $A^\N$ onto a subring of $A^\N$, and
the identities of \eqref{subsec_Witt-laws} tell us
that the ring structure induced on $A^\N$ via this
identification, is precisely the one given by the
addition law $S_A$ and multiplication law $P_A$, and
also the zero element and the unit are the required
ones, since $\bomega_A(\underline 0_A)=0\in A^\N$
and $\bomega_A(\underline 1_A)=1\in A^\N$. By the
same token, it is clear that $-\underline a$ is
computed by $I_A(\underline a)$ in the resulting ring.
\end{proof}

\begin{definition}\label{def_Witt-vectors}
Let $(A,\cT)$ be any topological ring. The topological ring
$$
(A^\N,S_A,P_A,\underline 0_A,\underline 1_A,\cT_{W(A)})
$$
is called {\em the ring of\/ Witt vectors associated with $A$},
and is denoted by $W(A,\cT)$. If no ambiguities are
likely to arise, we often omit reference to $\cT$, and
write simply $W(A)$.
\end{definition}

\sset\subsubsection{}
By construction, the ghost map is a continuous
ring homomorphism
$$
\bomega_A:W(A)\to A^\N
$$
(where $A^\N$ is endowed with termwise addition
and multiplication).
Henceforth, the addition and multiplication of elements
$\underline a,\underline b\in W(A)$ shall be denoted
simply in the usual way : $\underline a+\underline b$
and $\underline a\cdot\underline b$, and the neutral
element of addition and multiplication shall be denoted
respectively $0$ and $1$. The rule $A\mapsto W(A)$ is
a functor from the category of topological rings (and
continuous ring homomorphisms) to itself; indeed, if
$\psi:A\to B$ is any continuous ring homomorphism,
the map $\psi^\N$ yields a continuous ring homomorphism
$$
W(\psi):W(A)\to W(B)
\qquad
\underline a\mapsto(\psi(a_n)~|~n\in\N).
$$
and moreover one has the identity :
$$
\psi^\N\circ\bomega_A=\bomega_B\circ W(\psi).
$$

\begin{remark}\label{rem_non-unital-ring}
(i)\ \
More generally, if $A$ is not necessarily unital ring, the proof
of theorem \ref{th_Witt} shows that the datum
$(A^\N,S_A,P_A,\underline 0_A)$ is a non-unital ring $W(A)$, that
we shall also call the ring of Witt vectors associated with $A$.
Then clearly the rule $A\mapsto W(A)$ extends to an endofunctor
of the category of (associative, commutative) non-unital rings.

(ii)\ \
If $A$ is a unital ring, and $f:A\to B,g:A\to C$ two unital
$A$-algebras, the natural ring homomorphisms
$B\to B\otimes_AC\leftarrow C$ induce a map of $W(A)$-algebras
\set\begin{equation}\label{eq_hamas}
WB\otimes_{WA}WC\to W(B\otimes_AC).
\end{equation}
Namely, it is the unique $W(A)$-linear map such that
$\underline b\otimes\underline c\mapsto
(P_n(\underline b\otimes 1,1\otimes\underline c)~|~n\in\N)$
for every $\underline b\in W(B)$ and $\underline c\in W(C)$,
with $\underline b\otimes 1:=(b_n\otimes 1~|~n\in\N)$, and
likewise for $1\otimes\underline c$.

(iii)\ \
More generally, let $A$ be a unital ring; we define a
{\em non-unital right $A$-algebra} as the datum of a right
$A$-module $B$ and a (commutative, associative) non-unital
ring structure on $B$ such that $(bb')\cdot a=b\cdot(b'a)$
for every $a\in A$ and every $b,b'\in B$. Likewise we define
{\em non-unital left $A$-algebras}. To every left non-unital
$A$-algebra $B$, we may attach a unital $A$-algebra whose
underlying $A$-module is $A\oplus B$, whose structure map
$A\to A\oplus B$ is the natural inclusion, and such that
$$
(a,b)\cdot(a',b'):=(aa',ab'+a'b+bb')
$$
and likewise for right $A$-algebras : the details are left
to the reader.

(iv)\ \
With this terminology, we may then extend (ii) to the case
where $B$ and $C$ are respectively right and left non-unital
$A$-algebras : in this case, the homomorphism \eqref{eq_hamas}
of non-unital rings is obtained as follows. Recall that the
polynomial $P_n$ is bihomogeneous of degree $(p^n,p^n)$ for every
$n\in\N$ (remark \ref{rem_homogeneous-laws}(i)); hence there exists
$N\in\N$ such that $P_n=\sum_{i=1}^Nr_iP_{n,i}$, with $r_i\in\Z$,
and where each $P_{n,i}$ is a monomial of bidegree $(p^n,p^n)$,
{\em i.e.} there exists $M_i\in\N$ with
$P_{n,i}=\prod_{j=0}^{M_i}X_j^{d_j}Y_j^{e_j}$ and
$\sum_{j=0}^{M_i}d_jp^j=\sum_{j=0}^{M_i}e_jp^j=p^n$. Then we set
$Q_n(\underline b,\underline c):=
\sum_{i=1}^Nr_iP_{n,i}(\underline b,\underline c)$, with
$P_{n,i}(\underline b,\underline c):=
\prod_{j=0}^{M_i}b_j^{d_j}\otimes\prod_{j=0}^{M_i}c_j^{e_j}$ for
every $i=1,\dots,N$ and $\underline b,\underline c$ as in (ii);
lastly, \eqref{eq_hamas} is given by the rule :
$\underline b\otimes\underline c\mapsto
(Q_n(\underline b,\underline c)~|~n\in\N)$. In order to
check that \eqref{eq_hamas} is a map of non-unital rings,
endow $B':=A\oplus B$ and $C':=A\oplus C$ of the unital
$A$-algebra structures described in (iii), and notice that
\eqref{eq_hamas} is the restriction of the map of unital
$W(A)$-algebras $WB'\otimes_{WA}WC'\to W(B'\otimes_AC')$ of (ii).
\end{remark}

\sset\subsubsection{}\label{subsec_V-and-F}
Next, one defines two continuous mappings $W(A)\to W(A)$,
by the rule :
$$
\begin{aligned}
F_A(\underline a) :=\: & (F_n(a_0,\dots,a_{n+1})~|~n\in\N) \\
V_A(\underline a) :=\: & (0,a_0,a_1,\dots).
\end{aligned}
$$
$F_A$ and $V_A$ are often called, respectively, the {\em Frobenius\/}
and {\em Verschiebung\/} maps (\cite[Ch.IX, \S1, n.5]{BouAC}).
Consider also the mappings $f_A,v_A:A^\N\to A^\N$ such that
$$
f_A(\underline a):=(a_1,a_2,\dots)
\qquad
v_A(\underline a):=(0,pa_0,pa_1,\dots)
\qquad
\text{for every $\underline a:=(a_0,a_1,\dots)\in A^\N$}.
$$
Then $f_A$ is a continuous ring endomorphism of $A^\N$, and
$v_A$ is a continuous endomorphism of the additive group of
$A^\N$, and the basic identity characterizing $(F_n~|~n\in\N)$
can be written as :
\set\begin{equation}\label{eq_basic-F}
\bomega_A\circ F_A=f_A\circ\bomega_A.
\end{equation}
On the other hand, \eqref{eq_witt-pols} yields the identity :
\set\begin{equation}\label{eq_basic-V}
\bomega_A\circ V_A=v_A\circ\bomega_A.
\end{equation}
Let $\psi:A\to B$ be a continuous ring homomorphism;
directly from the definitions we obtain :
\set\begin{equation}\label{eq_functor_F-V}
W(\psi)\circ F_A=F_B\circ W(\psi)\qquad
W(\psi)\circ V_A=V_B\circ W(\psi).
\end{equation}

\begin{proposition}\label{prop_V_A-and-F_A}
With the notation of \eqref{subsec_V-and-F}, we have :
\begin{enumerate}
\item
The map $F_A$ is a continuous endomorphism of the ring $W(A)$,
and $V_A$ is a continuous endomorphism of the additive group
underlying $W(A)$.
\item
Moreover one has the identities :
\set\begin{equation}\label{eq_Verschiebung}
\begin{aligned}
F_A\circ V_A =\: & p\cdot\one_{W(A)} \\
V_A(\underline a\cdot F_A(\underline b)) =\: &
V_A(\underline a)\cdot\underline b
\quad\text{for every $\underline a,\underline b\in W(A)$}.
\end{aligned}
\end{equation}
\item
$F_A$ is a lifting of the Frobenius endomorphism of\/
$W(A)/pW(A)$, {\em i.e.} :
$$
F_A(\underline a)\equiv\underline a^p\pmod{pW(A)}
\qquad\text{for every $\underline a\in W(A)$.}
$$
\end{enumerate}
\end{proposition}
\begin{proof}(i): Arguing as in the proof of theorem
\ref{th_Witt}, we reduce to the case where $p$ is
a regular element in $A$ and the Frobenius endomorphism
of $A/pA$ lifts to a ring endomorphism $\phi:A\to A$,
in which case $\bomega_A$ is an injective ring homomorphism
(proposition \ref{prop_basic-cong}(i)). Then the assertion
follows immediately from \eqref{eq_basic-F} and
\eqref{eq_basic-V}. By the same token, we also get
the identities of (ii). Likewise, (iii) is reduced
to the identity
$$
f_A(\underline a)-\underline a^p\in p\cdot\Img\,\bomega_A
\qquad
\text{for every $\underline a:=(a_n~|~n\in\N)\in\Img\,\bomega_A$}
$$
(where here $\underline a^p=(a^p_n~|~n\in\N)$ denotes the
$p$-power map for the ring structure of $A^\N$). Now, we have
\set\begin{equation}\label{eq_many-congruences}
\phi(a_n)\equiv a_n^p\pmod{pA}
\qquad
\phi(a_n)\equiv a_{n+1}\pmod{p^{n+1}A}
\qquad
\text{for every $n\in\N$}
\end{equation}
by proposition \ref{prop_basic-cong}(ii), hence
$a_n^p-a_{n+1}\in pA$ for every $n\in\N$. Set
$b_n:=p^{-1}\cdot(a_n^p-a_{n+1})$; again by proposition
\ref{prop_basic-cong}(ii), it suffices to check that
$\phi(b_n)\equiv b_{n+1}\pmod{p^{n+1}A}$, that is
$$
\phi(a_n^p-a_{n+1})\equiv a_{n+1}^p-a_{n+2}
\pmod{p^{n+2}A}
\qquad
\text{for every $n\in\N$}.
$$
However, from \eqref{eq_many-congruences} and lemma
\ref{lem_basic-cong}(i) we deduce that
$\phi(a_n)^p\equiv a_{n+1}^p\pmod{p^{n+2}A}$, whence
the contention : details left to the reader.
\end{proof}

\sset\subsubsection{}\label{subsec_here-we-trunk}
Let $n\in\N$ be any integer; it follows easily from
\eqref{eq_Verschiebung} that
$$
V_n(A):=\Img\,V^n_A
$$
is an ideal of $W(A)$, and we define the ring of
{\em $n$-truncated Witt vectors\/} of $A$ as the quotient
$$
W_n(A):=W(A)/V_n(A)
$$
which we endow with the topology $\cT_{W_n(A)}$ induced by
the projection $\pi_n:W(A)\to W_n(A)$.

\begin{remark}\label{rem_truncated-Witt-non-unital}
(i)\ \
More generally, if $A$ is a not necessarily unital ring, then
$V_n(A)$ is a well-defined ideal of the non-unital ring $W(A)$
(see remark \ref{rem_non-unital-ring}(i)). Hence, the quotient
$W_n(A)$ is a well-defined non-unital ring in this case, for
every $n\in\N$.

(ii)\ \
Likewise, in the situation of remark \ref{rem_non-unital-ring}(iii),
we have a well-defined homomorphism of non-unital $A$-algebras
$W_nB\otimes_{W_nA}WC\to W_n(B\otimes_AC)$, for every $n\in\N$.
\end{remark}

\begin{lemma}\label{lem_Witt-truncate}
With the notation of \eqref{subsec_here-we-trunk}, we have :
\begin{enumerate}
\item
$\underline a=
(a_0,\dots,a_{m-1},0,\dots)+V^m_A\circ f^m_A(\underline a)$
for every $\underline a\in W(A)$ (where the addition is
taken in the additive group of\/ $W(A)$).
\item
The projection $W(A)\to A^n\ :\
(a_i~|~i\in\N)\mapsto(a_0,\dots,a_{n-1})$ factors through
$\pi_n$ and a bijection $W_n(A)\isom A^n$, and $\cT_{W_n(A)}$
corresponds to the product topology of $A^n$, under this
identification.
\end{enumerate}
\end{lemma}
\begin{proof}(i): Arguing as in the proof of theorem
\ref{th_Witt}, we may reduce to the case where $p$
is a regular element in $A$, in which case $\bomega_A$
is an injective ring homomorphism, and therefore it
suffices to show that
$$
\bomega_A(\underline a)=
\bomega_A(a_0,\dots,a_{m-1},0,\dots)+
\bomega_A\circ V^m_A\circ f^m_A(\underline a).
$$
The latter follows by a simple inspection : details left
to the reader. Assertion (ii) follows easily from (i).
\end{proof}

\begin{remark}\label{rem_Witt-limit}
(i)\ \
In view of lemma \ref{lem_Witt-truncate}(ii), we get
therefore an inverse system of topological rings
$(W_n(A)~|~n\in\N)$, whose limit is $W(A)$. Clearly
the ghost components descend to well-defined continuous
ring homomorphisms :
$$
\bar\bomega_m:W_n(A)\to A
\qquad\text{for every $m<n$}.
$$
Especially, the map $\bar\bomega_0:W_1(A)\to A$ is an
isomorphism of topological rings.

(ii)\ \
A direct inspection of the construction shows that the
endofunctor $W$ of the category of topological rings
commutes with all limits.

(iii)\ \
In the same vein, if the topological ring $A$ is the
limit of a system $(A_i~|~i\in I)$ of topological rings,
then $F_A$ is the limit of the corresponding system of
map $(F_{A_i}~|~i\in I)$, and likewise for $V_A$.

(iv)\ \
Let $A$ be any topological ring, $I$ any ideal of
the monoid $(A,\cdot)$, and $r\in\N$ any integer. In
light of remark \ref{rem_homogeneous-laws}(i), it is
easily seen that the sets
$$
W(I,r):=\{(a_n~|~n\in\N)\in W(A)~|~
        a_n\in I^{p^n}\!A\ \text{for every $n\leq r$}\}
\quad\text{and}\quad
W(I):=\bigcap_{n\in\N}W(I,n)
$$
are ideals of $W(A)$ : details left to the reader.
Moreover, we have the inclusions :
\set\begin{equation}\label{eq_general-inclusion}
W(I^c,r)\subset W(I,r)^c
\quad\text{and}\quad
W(I^c)\subset W(I)^c
\end{equation}
(where, for any topological space $X$ and every subset
$S\subset X$ we denote by $S^c$ the topological closure
of $S$ in $X$). Indeed, by definition we have
$W(I,r):=\prod_{n\in\N}J_n$, with $J_n:=I^{p^n}A$ for
every $n\leq r$ and $J_n:=A$ for $n>r$, so that
$W(I,r)^c=\prod_{n\in\N}J_n^c$, and it suffices to recall
that $(I^c)^{p^n}\subset(I^{p^n})^c$ for every $n\in\N$.
The same argument applies to $W(I)$.

Especially, if $I^{p^n}A$ is a closed ideal in the topology
of  $A$ for every $n\in\N$, then it follows that $W(I,r)$
is a closed ideal of $W(A)$, and the same holds also for
$W(I)$.

Furthermore, if $I,J$ are any two ideals of the monoid
$(A,\cdot)$, and $r\in\N$ is any integer, it follows easily
from remark \ref{rem_homogeneous-laws}(i) that
$$
W(I,r)\cdot W(J,r)\subset W(IJ,r)
\qquad\text{and}\qquad
W(I)\cdot W(J)\subset W(IJ).
$$
\end{remark}

\begin{proposition}\label{prop_int-with-nilker}
Let $A$ be any topological ring, $n\in\N$ any integer.
We have :
\begin{enumerate}
\item
The map (notation of remark {\em\ref{rem_Witt-limit}(i)})
$$
\pi_n:W_{n+1}(A)\to A^{n+1}
\qquad
\underline a\mapsto
(\bar\bomega_0(\underline a),\bar\bomega_1(\underline a),
\dots,\bar\bomega_n(\underline a))
$$
is an integral ring homomorphism with nilpotent kernel.
\item
$p^nA^{n+1}\subset\Img\,\pi_n\subset\{(x_0,\dots,x_n)\in A^{n+1}~|~
x_i\equiv x_0^{p^i}\pmod{pA}\ \ \text{for $i=1,\dots,n$}\}$.
\end{enumerate}
\end{proposition}
\begin{proof}(i): Denote by $B\subset A^{n+1}$ the integral
closure of the image of $\pi_n$. Clearly $B$
contains the $n+1$ canonical idempotents $e_0,\dots,e_n$
of $A^{n+1}$. Next, a simple inspection shows that
$$
\pi_n(ae_0)=(a,a^p,\cdots,a^{p^n})
\qquad
\text{for every $a\in A$}.
$$
It follows that
$a^{p^i}e_i=e_i\cdot\pi_n(ae_0)\in B$, and
then also $ae_i\in B$, for every $a\in A$ and every
$i\leq n$. This implies that $B=A^{n+1}$, {\em i.e.}
$\pi_n$ is integral, as asserted.
Next, the sequence of ideals $(V_i(A)~|~i\in\N)$
induces a descending filtration $(\bar V_i~|~i=0,\dots,n)$
on $W_{n+1}(A)$, and we notice:

\begin{claim}\label{cl_W_n-structure}
For every $i=0,\dots,n-1$, we have an isomorphism
of $W_{n+1}(A)$-modules
$$
A\isom\bar V_i/\bar V_{i+1}
\qquad
a\mapsto ae_i \pmod{\bar V_{i+1}}
$$
for the $W_{n+1}(A)$-module structure on $A$ induced by
$\bar\bomega_i$ (notation of \eqref{rem_Witt-limit}).
\end{claim}
\begin{pfclaim} It is clear that this map is
bijective. To show that it is also $W_{n+1}(A)$-linear,
notice that $P_j(X_0,\dots,X_j,0,\dots,0)=0$ for every
$j<i$ (remark \ref{rem_homogeneous-laws}(i)), so that
$$
\bomega_i(0,\dots,0,P_i(X_0,\dots,X_i,0,\dots,0,Y_i))=
\bomega_i(X_0,\dots,X_i)\cdot\bomega_i(0,\dots,0,Y_i)
$$
by \eqref{subsec_Witt-laws}. This translates as the
identity
$$
p^i\cdot P_i(X_0,\dots,X_i,0,\dots,0,Y_i)=
\bomega_i(X_0,\dots,X_i)\cdot p^iY_i
$$
whence $P_i(X_0,\dots,X_i,0,\dots,0,Y_i)=
\bomega_i(X_0,\dots,X_i)\cdot Y_i$, which implies the claim.
\end{pfclaim}

Set $I:=\Ker\,\pi_n$; from claim \ref{cl_W_n-structure}
we deduce that $I\cdot\bar V_i\subset\bar V_{i+1}$ for
every $i=0,\dots,n$, whence $I^{n+1}=0$ in $W_{n+1}(A)$,
which concludes the proof of (i).

(ii): The second inclusion follows by direct inspection
of the Witt polynomials. To show the first inclusion,
consider any sequence $(a_0,\dots,a_n)\in A^{n+1}$; we
need to find a solution $(x_0,\dots,x_n)\in W_{n+1}$
for the system of polynomial equations :
$$
x_0^{p^i}+px_1^{p^{i-1}}+\cdots+p^ix_i=p^na_i
\qquad
\text{for $i=0,\dots,n$}
$$
and we claim that there exists a solution that fulfills as
well the further condition : $x_i\in p^{n-i}A$ for $i=0,\dots,n$.
The proof proceeds by a simple induction on $i\leq n$, which
we leave to the reader.
\end{proof}

\begin{corollary}\label{cor_int-with-nilker}
Let $A$ be any topological ring such that $p^kA=0$ for
some $k\in\N$. Then $\Ker\,(\bar\omega_0:W_{n+1}A\to A)$
is nilpotent for every $n\in\N$.
\end{corollary}
\begin{proof} Let $\underline a\in W_{n+1}A$ be an element
such that $\bar\omega_0(\underline a)=0$, and define
$\pi_n$ as in proposition \ref{prop_int-with-nilker}(i);
in view of \eqref{eq_witt-pols}, we see that
$\pi_n(\underline a)=(0,pb_1,\dots,pb_n)$ for some
$b_1,\dots,b_n\in A$; hence $\underline a^k\in\Ker\,\pi_n$,
and it suffices to invoke proposition
\ref{prop_int-with-nilker}(i) to conclude.
\end{proof}

\begin{lemma}\label{lem_Witt-limit}
Let $(A,\cT)$ be any topological ring, with separated
completion $(A^\wedge,\cT^\wedge)$.
\begin{enumerate}
\item
If the topology $\cT$  is discrete, the topology
$\cT_{W(A)}$ agrees with the linear topology defined
by the filtration $(V_m(A)~|~m\in\N)$.
\item
If the topology $\cT$ is separated (resp. and complete),
then the same holds for $\cT_{W(A)}$.
\item
If $\cT$ is a linear topology, then the same holds
for $\cT_{W(A)}$.
\item
The natural map $W(A,\cT)\to W(A^\wedge,\cT^\wedge)$
factors uniquely through a natural isomorphism of topological
rings
$$
W(A,\cT)^\wedge\isom W(A^\wedge,\cT^\wedge)
$$
(where $W(A,\cT)^\wedge$ denotes the separated completion
of\/ $W(A,\cT)$).
\end{enumerate}
\end{lemma}
\begin{proof} (i) follows easily from lemma
\ref{lem_Witt-truncate}(ii) : details left to the reader.

(ii): It is clear from the definitions that if $A$ is
separated, then the same holds for $\cT_{W(A)}$. Next,
suppose that $A$ is complete and separated, and for
every $n\in\N$ consider the short exact sequence of
topological groups
\set\begin{equation}\label{eq_reduce-to-W_n}
0\to V_n(A)/V_{n+1}(A)\to W_{n+1}(A)\to W_n(A)\to 0
\end{equation}
where $V_n(A)/V_{n+1}(A)$ is endowed with the topology
induced by the inclusion into $W_{n+1}(A)$. Taking into
account claim \ref{cl_W_n-structure} and proposition
\ref{prop_replaces-Mat-Th.8.1}(i,v), a simple induction
then shows that $W_n(A)$ is complete and separated for
every $n\in\N$. Then the assertion follows from corollary
\ref{cor_limits-and-complete}(i) and remark
\ref{rem_Witt-limit}(i).

(iii): Let $(I_\lambda~|~\lambda\in\Lambda)$ be a
fundamental system of open ideals of $A$, and for every
$\lambda\in\Lambda$ and every $r\in\N$, denote by
$\pi_{\lambda,r}:W(A)\to W_r(A/I_\lambda)$ the projection;
then it is easily seen that the family of ideals
$(\Ker\,\pi_{\lambda,r}~|~\lambda\in\Lambda,\ r\in\N)$
is a fundamental system of open neighborhoods of $0\in W(A)$.
The assertion is an immediate consequence.

(iv) is similar to (ii) : in light of corollary
\ref{cor_limits-and-complete}(ii), we are reduced
to checking that $W_n(A^\wedge)$ is the separated
completion of $W_n(A)$. The latter is established by
induction on $n\in\N$, by means of the exact sequences
\eqref{eq_reduce-to-W_n} and proposition
\ref{prop_replaces-Mat-Th.8.1}(i,v).
\end{proof}

\sset\subsubsection{}\label{subsec_Teich}
The projection $\bomega_0$ admits a continuous set-theoretic
splitting, the {\em Teichm{\"u}ller mapping}
$$
\tau_A:A\to W(A) \qquad a\mapsto(a,0,0,\dots).
$$
For every $a\in A$, we call $\tau_A(a)$ the
{\em Teichm\"uller representative\/} of $a$.
Clearly, for any continuous ring homomorphism
$\psi:A\to B$ we have :
\set\begin{equation}\label{eq_Teich-functorial}
\tau_B\circ\psi=W(\psi)\circ\tau_A.
\end{equation}

\begin{proposition}\label{prop_Teich-series}
With the notation of \eqref{subsec_Teich}, we have :
\begin{enumerate}
\item
$\tau_A(a)\cdot\underline b=(a^{p^n}b_n~|~n\in\N)$ for
every $a\in A$ and $\underline b:=(b_n~|~n\in\N)\in W(A)$.
\item
Moreover, $\tau_A$ is a multiplicative map, {\em i.e.}
the following identity holds :
$$
\tau_A(a)\cdot\tau_A(b)=\tau_A(a\cdot b)\qquad
\text{for every $a,b\in A$}
$$
\item
$\underline a=\sum^\infty_{n=0}V_A^n(\tau_A(a_n))$
for every $\underline a:=(a_n~|~n\in\N)\in W(A)$,
where the convergence of the series is relative to
the topology $\cT_{W(A)}$.
\end{enumerate}
\end{proposition}
\begin{proof}(i): Since $P_n$ is bihomogeneous of degree
$(p^n,p^n)$ (remark \ref{rem_homogeneous-laws}(i)), we have
$$
P_n(X_0,0,\dots,0,Y_0,\dots,Y_n)=
X_0^{p^n}\cdot P_n(1,0,\dots,0,Y_0,\dots,Y_n)=
X_0^{p^n}\cdot Y_n
$$
where the second identity follows after recalling that
$\underline 1_A:=(1,0,\dots)$ is the unit of $W(A)$. The
assertion follows immediately. (ii) is a special case of
(i), and (iii) follows easily from lemma
\ref{lem_Witt-truncate}(i).
\end{proof}

\begin{example}\label{ex_how-to-write-p-powers}
(i)\ \ 
For any $n\in\N$, let $(a_0,\dots,a_n)\in W_{n+1}(\Z)$ be
the sequence that represents the element $p^n$. Since
$\bar\omega_i$ is a ring homomorphism, we have
$$
p^n=a_0^{p^i}+pa_1^{p^{i-1}}+\cdots+p^ia_i
\qquad
\text{for every $i=0,\dots,n$}.
$$
A simple induction then shows that
$$
a_i\equiv p^{n-i}\pmod{p^{n-i+1}}
\qquad
\text{for every $i=0,\dots,n$}.
$$

(ii)\ \
Let $A$ be any topological ring, $n,i\in\N$ any integers;
set
$$
k:=\min(i,n)
\qquad\text{and}\qquad
\bar V_{n,j}(A):=V_j(A)/V_{n+1}(A)
\qquad
\text{for every $j\leq n$}.
$$
Let also $\bar V:W_n(A)\to W_{n+1}(A)$ be the $\Z$-linear
map induced by $V$. We claim that
$$
\sum_{j=0}^kp^{i-j}\bar V_{n,j}(A)=J_{i,n}:=
\{(a_0,\dots,a_n)\in W_{n+1}(A)~|~
a_j\in p^{i-j}A\ \text{for every}\ j\leq k\}.
$$
For the proof, we argue by induction on $n$ and $i$, and
notice first that $J_{i,n}$ is an ideal of $W_{n+1}(A)$ for
every $n$ and $i$, due to remark \ref{rem_homogeneous-laws}(i).
If $n=0$ or $i=0$, there is nothing to prove. Suppose then
that $n,i>0$ and that the assertion is known for all
strictly smaller values of $n$ and $i$; from (i) we see
that $p^iW_{n+1}(A)\subset J_{i,n}$, and the inductive
assumption shows that
$$
\sum_{j=1}^kp^{i-j}\bar V_j(A)=
\bar V\Bigl(\sum_{j=0}^{k-1}p^{i-j-1}\bar V_{n-1,j}(A)\Bigr)=
\bar V(J_{i-1,n-1})\subset J_{i,n}
$$
whence $\sum_{j=0}^kp^{i-j}\bar V_{n,j}(A)\subset J_{i,n}$.
Conversely, say that $\underline a:=(a_0,\dots,a_n)\in J_{i,n}$,
so that $a_0=p^ib_0$ for some $b_0\in A$; due to lemma
\ref{lem_Witt-truncate}(i) and proposition
\ref{prop_Teich-series}(i) we may write
$$
\underline a=\tau(a_0)+\bar V(a_1,\dots,a_n)=
\tau(b_0)\cdot(p^i,0,\dots,0)+\bar V(a_1,\dots,a_n)
$$
and notice that $\bar V(a_1,\dots,a_n)\in
\bar V(J_{i-1,n-1})\subset\sum_{j=1}^kp^{i-j}\bar V_{j,n}(A)$.
Also, from (i) we see that
$(p^i,0,\dots,0)=p^i+\bar V(b_1,\dots,b_n)$, where again
$(b_1,\dots,b_n)\in J_{i-1,n-1}$. Summing up, we get
$J_{i,n}\subset\sum_{j=0}^kp^{i-j}\bar V_{j,n}(A)$, as required.

(iii)\ \
Especially, (ii) implies that for every ring $A$ we have
$$
p^iW_{n+1}(A)\subset J_{i,n}\subset p^{i-n}W_{n+1}(A)
\qquad
\text{for every $i,n\in\N$ with $i\geq n$}.
$$
If $\cT_p$ denotes the $p$-adic topology on $A$, it
follows that the topology of $W_{n+1}(A,\cT_p)$ agrees
with the $p$-adic topology of the underlying ring
$W_{n+1}(A)$, for every $n\in\N$.
\end{example}

\begin{example}\label{ex_localize-Witt}
(i)\ \ 
Let $S\subset\Z$ be any subset, and $A$ any ring (which
we may regard as a discrete topological ring). Then the
localization map $j:A\to S^{-1}A$ induces a ring isomorphism
$$
S^{-1}W_{n+1}(A)\isom W_{n+1}(S^{-1}A)
\qquad
\text{for every $n\in\N$}.
$$
Indeed, let us show first that multiplication by every
$s\in S$ is a bijection on $W_{n+1}(S^{-1}A)$. To this aim, we
consider the descending filtration $(\bar V_iA~|~i=0,\dots,n)$
as in the proof of proposition \ref{prop_int-with-nilker}(i);
we are then easily reduced to checking that multiplication
by $s$ is bijective on $\bar V_i(A)/\bar V_{i+1}(A)$ for
$i=0,\dots,n-1$. But according to claim \ref{cl_W_n-structure},
the latter quotient is isomorphic to $S^{-1}A$, for the
$W_{n+1}(S^{-1}A)$-module structure induced on $S^{-1}A$ by
$\bar\omega_i$, whence the contention. Thus, $W_{n+1}(j)$
extends to a well defined ring homomorphism
$\lambda:S^{-1}W_{n+1}A\to W_{n+1}(S^{-1}A)$. To see that the
latter is an isomorphism, we notice that $\lambda$ is a
homomorphism of filtered rings for the filtrations
$S^{-1}\bar V_\bullet A$ on $S^{-1}W_{n+1}A$ and
$\bar V_\bullet(S^{-1}A)$ on $W_{n+1}(S^{-1}A)$; then it
suffices to check that the associated graded ring
homomorphism $\gr_\bullet\lambda$ is bijective, and the
latter assertion follows again from claim \ref{cl_W_n-structure}.

(ii)\ \
If $p$ is invertible in $A$, a simple inspection shows
that ghost map $\bomega_A:W(A)\to A^\N$ is a ring homomorphism,
and likewise, we have $W_{n+1}A\isom A^{n+1}$ in this case,
for every $n\in\N$. Combining with (i), we deduce that
for every ring $A$, the localization map $A\to A[p^{-1}]$
induces a ring isomorphism
$$
W_{n+1}(A)[p^{-1}]\isom (A[p^{-1}])^{n+1}
\qquad
\text{for every $n\in\N$}.
$$

(iii)\ \
In the same vein, let $S\subset A$ be any multiplicative
subset, and denote by $\tau_A(S)\subset W_{n+1}A$ the image
of the subset $\{\tau_A(s)~|~s\in S\}\subset W(A)$ (notation
of \eqref{subsec_Teich}). Then the localization map
$j:A\to S^{-1}A$ induces a ring isomorphism
$$
\lambda:\tau_A(S)^{-1}W_{n+1}(A)\isom W_{n+1}(S^{-1}A)
\qquad
\text{for every $n\in\N$}.
$$
Indeed, from proposition \ref{prop_Teich-series}(ii)
we see that the image of $\tau_A(S)$ is invertible in
$W_{n+1}(S^{-1}A)$, so $\lambda$ is well defined. Next,
proposition \ref{prop_Teich-series}(i) easily implies
that $\lambda$ is surjective. Lastly, the kernel of
$\lambda$ is $\tau_A(S)^{-1}\Ker\,W_{n+1}(j)$; but clearly
$\Ker\,W_{n+1}(j)$ is the ideal of all sequences
$\underline a:=(a_0,\dots,a_n)$ such that for every
$i=0,\dots,n$ there exists $s_i\in S$ with $s_ia_i=0$.
Then $\tau_A(s_0\cdots s_n)\cdot\underline a=0$ in
$W_{n+1}A$; this shows that $\Ker\,\lambda=0$, and
concludes the proof.
\end{example}

\sset\subsubsection{}\label{subsec_filtrations-agree}
Let now $A$ be a topological $\F_p$-algebra. In this case
\eqref{eq_congr-for-F_n} yields the identity :
\set\begin{equation}\label{eq_simple-F}
F_A(\underline a)=(a^p_n~|~n\in\N)
\qquad\text{for every $\underline a:=(a_n~|~n\in\N)\in W(A)$}.
\end{equation}
As an immediate consequence we get :
\set\begin{equation}\label{eq_F-V-commute}
p\cdot\underline a=V_A\circ F_A(\underline a)=F_A\circ V_A(\underline a)=
(0,a_0^p,a_1^p,\dots)
\qquad\text{for every $\underline a\in W(A)$}.
\end{equation}
Especially, we have
\set\begin{equation}\label{eq_this-is-p}
p=(0,1,0,\dots)
\qquad
\text{in $W(A)$}.
\end{equation}
Moreover, let $\Phi_A:A\to A$ be the Frobenius endomorphism;
\eqref{eq_simple-F} also implies the identity :
\set\begin{equation}\label{eq_Teich-Frob}
F_A\circ\tau_A=\tau_A\circ\Phi_A.
\end{equation}

\begin{proposition}\label{prop_Witt-is-complete}
Let $A$ be any $\F_p$-algebra. We have :
\begin{enumerate}
\item
The $p$-adic topology on $W(A)$ agrees with the $V_1(A)$-adic
topology, and both are finer than the topology $\cT'$ defined
by the filtration $(V_n(A)~|~n\in\N)$.
\item
Moreover, $W(A)$ is also complete and separated for the
$p$-adic topology.
\item
$W(A)^\times=\bomega_0^{-1}(A^\times)$.
\item
$V^i(\tau_A(b))\cdot V^j(\tau_A(b))=
V^{i+j}(\tau_A(a^{p^j}b^{p^i}))$
\ \ for every $a,b\in A$ and every $i,j\in\N$.
\end{enumerate}
\end{proposition}
\begin{proof}(i): It is clear from \eqref{eq_F-V-commute}
that the $p$-adic topology on $W(A)$ is finer than both
$\cT'$ and the $V_1$-adic topology. Conversely, for every
$\underline a,\underline b\in W(A)$ we have
$$
V_A(\underline a)\cdot V_A(\underline b)=
V_A(\underline a\cdot F_AV_A(\underline b))=
p\cdot V_A(\underline a\cdot\underline b)
$$
whence $V_1^2\subset pA$, which says that the $V_1$-adic
topology is finer than the $p$-adic topology.

(ii): Since $W(A)$ is complete and separated for the
linear topology $\cT'$, the assertion follows from lemma
\ref{lem_fontaine}.

(iii): Clearly $W(A)^\times\subset\bomega^{-1}_0(A^\times)$.
For the converse, notice that the ideal $V_1(A)$ is
contained in the Jacobson radical ideal of $W(A)$, due
to (i),(ii) and remark \ref{rem_someth-on-bdd-in-Z-lin}(v).
Therefore, the induced map
$\Spec\,\bomega_0:\Spec\,A\to\Spec\,W(A)$ restricts to
a bijection on the sets of maximal ideals of $A$ and
$W(A)$, and the assertion follows easily.

(iv): Without loss of generality, we may assume that
$i\leq j$; then we compute
$$
\begin{aligned}
V^i(\tau_A(b))\cdot V^j(\tau_A(b))=\, &
V^i(\tau_A(a)\cdot F^iV^j(\tau_A(b))) & &
\quad\text{(by \eqref{eq_Verschiebung})} \\
=\, & V^i(\tau_A(a)\cdot p^i\cdot V^{j-i}(\tau_A(b)))  & &
\quad\text{(again by \eqref{eq_Verschiebung})} \\
=\, & p^i\cdot V^i(\tau_A(a)\cdot V^{j-i}(\tau_A(b))) & &
\quad\text{(by proposition \ref{prop_V_A-and-F_A}(i))}\\
=\, & p^i\cdot V^i(V^{j-i}(\tau_A(a^{p^{j-i}}b))) & &
\quad\text{(by proposition \ref{prop_Teich-series}(i))}\\
=\, & p^i\cdot V^j(\tau_A(a^{p^{j-i}}b)) & &
\quad\text{(by \eqref{eq_F-V-commute})} \\
=\, & V^{i+j}(\tau_A(a^{p^j}b^{p^i}))
\end{aligned}
$$
as stated.
\end{proof}

\sset\subsubsection{}\label{subsec_perfect-case}
Suppose additionally that $A$ is {\em perfect}, {\em i.e.}
that $\Phi_A$ is an automorphism of the topological ring
$A$ (especially, it is a homeomorphism). It follows from
\eqref{eq_simple-F} that $F_A$ is an automorphism of $W(A)$.
Hence, in view of \eqref{eq_Verschiebung} and
\eqref{eq_F-V-commute} :
$$
p^n\cdot W(A)=V_n(A)
\qquad\text{for every $n\in\N$}.
$$
So, the $p$-adic filtration coincides with the filtration
$(V_n(A)~|~n\in\N)$ and with the $V_1(A)$-adic filtration.
Especially, the $0$-th ghost component descends to an
isomorphism
$$
\bar\bomega_0:W(A)/pW(A)\isom A.
$$
Also, the identity of proposition \ref{prop_Teich-series}(iii)
can be written in the form :
\set\begin{equation}\label{eq_new-form}
\underline a=\sum^\infty_{n=0}p^n\cdot\tau_A(a_n^{p^{-n}})
\qquad\text{for every $\underline a:=(a_n~|~n\in\N)\in W(A)$}.
\end{equation}

\begin{proposition}\label{prop_reduced-Witt}
Let $A$ be a reduced $\F_p$-algebra,
$\underline a:=(a_n~|~n\in\N)\in W(A)$ any element.
\begin{enumerate}
\item
Consider the following conditions :
\begin{enumerate}
\item
The ideal $J:=\sum_{n\in\N}a_nA$ is not contained in
any minimal prime ideal of $A$.
\item
$\underline a$ is a regular element of\/ $W(A)$.
\end{enumerate}
Then {\em(a)$\Rightarrow$(b)}, and if $A$ has only
finitely many minimal prime ideals, {\em(b)$\Rightarrow$(a)}.
\item
$W(A)$ is reduced, and if $A$ is a domain, the
same holds for $W(A)$.
\item
If $A$ is a perfect field, then $W(A)$ is a complete
discrete valuation ring of mixed characteristic $(0,p)$,
maximal ideal $pW(A)$, and residue field $A$.
\item
If $A\subset B$ is an inclusion of perfect $\F_p$-algebras,
then $W(A)\subset W(B)$, and
$$
p^nW(A)=W(A)\cap p^nW(B)
\qquad
\text{for every $n\in\N$}.
$$
\item
If\/ $A$ is perfect, and\ \ $\sum_{i=0}^ka_iA=A$\ \ for some
$k\in\N$, we have
$$
\underline aW(A)\cap p^nW(A)=
p^{n-k}\cdot(\underline aW(A)\cap p^kW(A))
\qquad
\text{for every $n\in\N$}
$$
and moreover, the ideal $\underline aW(A)$ is a closed
subset for the $p$-adic topology on $W(A)$.
\end{enumerate}
\end{proposition}
\begin{proof} Notice that the assertions are independent of
the topology of $A$ and $B$, hence we may assume that both
these topologies are discrete.

(iii): We know already that $W(A)$ is $p$-adically
complete, by proposition \ref{prop_Witt-is-complete}(ii).
Also, it follows easily from \eqref{eq_F-V-commute} that
$p$ is regular in $W(A)$, and taking into account
proposition \ref{prop_Witt-is-complete}(iii) we see that
every element of $W(A)$ can be written uniquely in the 
form $p^n\cdot\underline u$, for some $n\in\N$ and an
invertible element $\underline u\in W(A)$. The assertion
follows immediately.

(ii): Suppose first that $A$ is a domain, and we prove
that the same holds for $W(A)$. To this aim, let $B$
be any perfect field containing $A$; since $W(A)$ is
a subring of $W(B)$, it suffices to check that $W(B)$
is a domain, which is clear from (iii).

Next, let $\Min\,A$ be the set of all minimal prime
ideals of $A$; if $A$ is reduced, the natural map
$j:A\to B:=\prod_{\fp\in\Min\,A}A_\fp$ is injective,
and therefore $W(A)$ is a subring of
$W(B)=\prod_{\fp\in\Min\,A}W(A_\fp)$.
To prove that $W(A)$ is reduced, it then suffices to
show that the same holds for every factor $W(A_\fp)$,
but this is already known, since each such $A_\fp$ is
a field.

(i): Suppose that (a) holds; arguing as in the proof
of (ii), we reduce to checking that the image of
$\underline a$ is regular in $W(A_\fp)$, for every
$\fp\in\Min\,A$. But under condition (a), it is clear
that the image of $\underline a$ does not vanish in
$W(A_\fp)$ for any such $\fp$. Thus, in this case (ii)
implies that (b) holds. Conversely, suppose that
$\Min\,A=\{\fp_1,\dots,\fp_r\}$ is a finite set, and
say that $J\subset\fp_1$. Pick any
$f\in\bigcap_{i=2}^r\fp_i\setminus\fp_1$; since $j$
is injective, it is easily seen that $\Ann_A(f)=\fp_1$
(details left to the reader). By proposition
\ref{prop_Teich-series}(i), we then see that
$\underline a\cdot\tau_A(f)=0$, which contradicts (b).

(iv): It is clear that $W(A)\subset W(B)$, and
also that $V_n(A)=W(A)\cap V_n(B)$, so the assertion
follows from \eqref{subsec_perfect-case}.

(v): Say that $x\underline a=p^ny$ for some $x,y\in W(A)$;
by (i), $p$ is regular in $W(A)$, so for the proof of
the stated equality it suffices to show that
$x\in p^{n-k}W(A)$. However, we have already noticed
that the map $j$ is injective, and therefore
$p^kW(A)=W(A)\cap p^kW(B)$, by (iv).
Hence, it suffices to show that the image of $x$ lies
in $p^{n-k}W(A_\fp)$ for every such $\fp$. Moreover,
$A_\fp$ is a perfect field for every $\fp\in\Min\,A$,
and $\sum_{i=0}^ka_i\cdot A_\fp=A_\fp$, so we may replace
from start $A$ by any $A_\fp$, in which case $W(A)$ is
a discrete valuation ring with $p$ as uniformizer, by
(iii), and notice that by assumption the elements
$a_0,\dots,a_k$ cannot all vanish, therefore
$$
\underline a\notin p^{k+1}W(A)
\qquad\text{and}\qquad
x\cdot\underline a\in p^nW(A)
$$
whence the assertion : details left to the reader.
Lastly, the foregoing implies that the $p$-adic topology
on $\underline aW(A)$ agrees with the topology induced
from the $p$-adic topology on $W(A)$. On the other hand,
(i) says that $\underline a$ is regular on $W(A)$, hence
scalar multiplication by $a$ induces a bijection
$W(A)\isom\underline aW(A)$ that identifies the $p$-adic
topology of $W(A)$ with the $p$-adic topology of
$\underline aW(A)$.
Since $W(A)$ is $p$-adically complete and separated
(proposition \ref{prop_Witt-is-complete}(ii)), it
follows that the same holds for $\underline aW(A)$.
Summing up, $\underline aW(A)$ is complete and separated
for the topology induced from the $p$-adic topology of
$W(A)$, and consequently it is a closed subset in $W(A)$,
for this topology.
\end{proof}

\begin{example}\label{ex_discrete-Witt}
Let be $(A,\cT)$ a topological $\F_p$-algebra.

(i)\ \
If $A$ is perfect and $\cT$ is the discrete topology,
then it follows from \eqref{subsec_perfect-case} and
proposition \ref{prop_Witt-is-complete}(i) that
$\cT_{W(A)}$ agrees with the $p$-adic topology.

(ii)\ \
By inspecting the proof of theorem
\ref{th_complete-top-grps}(i), it is easily seen that if
$(A,\cT)$ is perfect, the same holds for its separated
completion.

(iii)\ \
For instance, we have a natural isomorphism of topological
rings :
$$
\Z_p\isom W(\F_p)
$$
(where $\F_p$ is endowed with the discrete topology, and
$\Z_p$ with the $p$-adic topology). Indeed, (i) says that
the topology of $W(\F_p)$ agrees with the $p$-adic topology,
hence the unique ring homomorphism $\Z\to W(\F_p)$ extends
to a unique continuous map $\psi:\Z_p\to W(\F_p)$. Since $p$
is a regular element of $W(\F_p)$ (proposition
\ref{prop_reduced-Witt}(i)), the map is injective.
Moreover, if we endow both rings with their $p$-adic
filtrations, the induced map
$\gr_\bullet\psi:\gr_\bullet\Z_p\to\gr_\bullet W(\F_p)$
on associated graded rings is easily seen to be an
isomorphism. Then the assertion follows from
\cite[Ch.III, \S2, n.8, Cor.2]{BouAC}.

(iv)\ \
It follows from (iii) that the structure map
$\beta:\F_p\to A$ induces a natural structure
$W(\beta)$ of topological $\Z_p$-algebra on $W(A)$.

(v)\ \
Let $A:=\F_p[T^{1/p^\infty}]$, the universal perfect
$\F_p$-algebra in one generator, and endow $A$ with
its discrete topology. Then $W(A)$ is the $p$-adic
completion $\Z_p\{T^{1/p^\infty}\}$ of $\Z_p[T^{1/p^\infty}]$,
together with its $p$-adic topology. The latter can be
described as the ring of all power series
\set\begin{equation}\label{eq_fract-powerseries}
\sum_{n\in\N}a_nT^{\lambda_n}
\end{equation}
where $(\lambda_n~|~n\in\N)$ is any sequence of non-negative
elements of $\Z[1/p]$, and $(a_n~|~n\in\N)$ is any sequence
of elements of $\Z_p$ that converges $p$-adically to $0$.
For the proof, we argue as in (iii): we have a unique map
of $\Z_p$-algebras $\psi:\Z_p[T^{1/p^\infty}]\to W(A)$ such
that $\psi(T^\lambda):=\tau_A(T^\lambda)$ for every $\lambda\geq 0$
in $\Z[1/p]$ (for the $\Z_p$-algebra structure on $W(A)$
provided by (iv)). Then $\psi$ extends by continuity
to a unique map of topological $\Z_p$-algebras
$\psi^\wedge:\Z_p\{T^{1/p^\infty}\}\to W(A)$, and if
we endow both rings with their $p$-adic filtrations,
it is easily seen that the resulting map of associated
graded rings $\gr^\bullet\psi^\wedge$ is an isomorphism,
so the same holds for $\psi^\wedge$, again by
\cite[Ch.III, \S2, n.8, Cor.3]{BouAC}.
\end{example}

\sset\subsubsection{}\label{subsec_liftings}
In the situation of lemma \ref{lem_basic-cong}, endow $A$
(resp. $A/J_n$, for every $n\in\N$) with the linear topology
$\cT$ defined by the system of ideals $(J_n~|~n\in\N)$
(resp. with the discrete topology) and let $\pi_n:A\to A/J_n$
and $\pi'_n:A/J_n\to A/J_0$ be the natural projections, for
every $n\in\N$. Let also $R$ be a topological $\F_p$-algebra,
and $\bar\phi:R\to A/J_0$ a continuous ring homomorphism. We
say that a continuous mapping $\phi_n:R\to A/J_n$ is a
{\em lifting of\/ $\bar\phi$\/} if we have :
$$
\phi_n(x^p)=\phi_n(x)^p
\qquad\text{for every $x\in R$, and}\qquad
\pi'_n\circ\phi_n=\bar\phi.
$$

\begin{lemma}\label{lem_liftings}
With the notation of \eqref{subsec_liftings}, the
following holds :
\begin{enumerate}
\item
If $\phi_n$ and $\phi'_n$ are two liftings of\/ $\bar\phi$,
then $\phi_n$ and $\phi'_n$ agree on $R^{p^n}$.
\item
If $R$ is a perfect topological $\F_p$-algebra, there exists
a unique lifting $\phi_n:R\to A/J_n$ of\/ $\bar\phi$. We have
$\phi_n(1)=1$ and $\phi_n(xy)=\phi_n(x)\cdot\phi_n(y)$ for
every $x,y\in R$.
\item
Let $R$ be as in {\em(ii)}, and suppose that $(A,\cT)$ is complete
and separated. Then :
\begin{enumerate}
\item
There exists a unique continuous mapping $\phi:R\to A$
such that
$$
\pi_0\circ\phi=\bar\phi
\qquad\text{and}\qquad
\phi(x^p)=\phi(x)^p
\quad
\text{for every $x\in R$}.
$$
\item
Moreover, $\phi(1)=1$ and $\phi(xy)=\phi(x)\cdot\phi(y)$
for every $x,y\in R$.
\item
If furthermore, $A$ is an $\F_p$-algebra, then
$\phi$ is a homomorphism of topological rings.
\end{enumerate}
\end{enumerate}
\end{lemma}
\begin{proof} (i) is an easy consequence of lemma
\ref{lem_basic-cong}(i).

(ii): The uniqueness follows from (i). For the existence,
pick any mapping $\sigma:A/J_0\to A/J_n$ such that
$\pi'_n\circ\sigma=\one_{A/J_0}$, and notice that the
mapping $\psi:=\sigma\circ\bar\phi:R\to A/J_n$ is
continuous, and $\pi'_n\circ\psi=\bar\phi$. We let :
$$
\phi_n(x):=\psi(x^{p^{-n}})^{p^n}
\qquad
\text{for every $x\in R$}.
$$
Using lemma \ref{lem_basic-cong}(i) one verifies easily that
$\phi_n$ does not depend on the choice of $\psi$. Especially,
define $\psi':R\to A/J_n$ by the rule : $x\mapsto\psi(x^{p^{-1}})^p$
for every $x\in R$. Clearly $\pi'_n\circ\psi'=\bar\phi$ as well,
hence :
$$
\phi_n(x^p)=\psi'(x^{p^{1-n}})^{p^n}=\psi(x^{p^{-n}})^{p^{n+1}}=\phi_n(x)^p
\qquad
\text{for every $x\in R$}
$$
as claimed. If we choose $\sigma$ so that $\sigma(1)=1$, we obtain
$\phi_n(1)=1$. Finally, since $\bar\phi$ is a ring homomorphism
we have $\psi(x)\cdot\psi(y)\equiv\psi(xy)\pmod{J_0}$ for every
$x,y\in R$. Hence
$$
\psi(x)^{p^n}\cdot\psi(y)^{p^n}\equiv\psi(xy)^{p^n}\pmod{J_n}
$$
(again by lemma \ref{lem_basic-cong}(i)) and finally
$\phi_n(x^{p^n})\cdot\phi_n(y^{p^n})=\phi_n((xy)^{p^n})$,
which implies the last stated identity, since $R$ is perfect.

(iii): The existence and uniqueness of $\phi$ follow from (ii).
It remains only to show that $\phi(x)+\phi(y)=\phi(x+y)$ in
case $A$ is an $\F_p$-algebra. It suffices to check the latter
identity on the projections onto $A/J_n$, for every $n\in\N$,
in which case one argues analogously to the foregoing proof
of the multiplicative property for $\phi_n$ : the details shall
be left to the reader.
\end{proof}

\begin{proposition}\label{prop_lift-Witt}
In the situation of \eqref{subsec_liftings}, suppose that
$R$ is a perfect topological $\F_p$-algebra and $(A,\cT)$
is complete and separated. Then we have :
\begin{enumerate}
\item
For every $n\in\N$ there exists a unique ring homomorphism
$v_n:W_{n+1}(A/J_0)\to A/J_n$ such that the following diagram
commutes :
$$
\xymatrix{
W_{n+1}(A) \ar[r]^-{\bomega_n} \ar[d]_{W_{n+1}(\pi_0)} &
A \ar[d]^{\pi_n} \\
W_{n+1}(A/J_0) \ar[r]^-{v_n} & A/J_n.
}$$
\item
Let $u_n:=v_n\circ W_{n+1}(\bar\phi)\circ\bar F{}_{\!\!R}^{-n}$,
where $\bar F_{\!R}:W_{n+1}(R)\to W_{n+1}(R)$ is induced by the
Frobenius automorphism of\/ $W(R)$. Then, for every $n\in\N$ the
following diagram commutes :
$$
\xymatrix{
W_{n+2}(R) \ar[r]^-{u_{n+1}} \ar[d] & A/J_{n+1} \ar[d]^\theta \\
W_{n+1}(R) \ar[r]^-{u_n} & A/J_n
}$$
where the vertical maps are the natural projections.
\item
There exists a unique ring homomorphism $u$ such that
the following diagram commutes:
$$
\xymatrix{
W(R) \ar[rr]^-u \ar[d]_{\bomega_0} & & A \ar[d]^{\pi_0} \\
R \ar[rr]^-{\bar\phi} & & A/J_0.
}$$
Furthermore, $u$ is continuous for the topology
$\cT_{W(R)}$, and we have :
\set\begin{equation}\label{eq_explicit}
u(\underline a)=\sum_{n\in\N}p^n\cdot\phi(a_n^{p^{-n}})
\qquad
\text{for every $\underline a:=(a_n~|~n\in\N)\in W(R)$}
\end{equation}
with $\phi:R\to A$ the unique continuous mapping
characterized as in lemma {\em\ref{lem_liftings}(iii)}.
\end{enumerate}
\end{proposition}
\begin{proof} (i): We have to check that
$\bomega_n(a_0,\dots,a_n)\in J_n$ whenever
$a_0,\dots,a_n\in J_0$, which is clear from lemma
\ref{lem_basic-cong}(ii).

(ii): Since $\bar F_{\!R}$ is an automorphism, it boils
down to verifying :
\begin{claim}
$\theta\circ v_{n+1}(\bar\phi(a_0),\dots,\bar\phi(a_{n+1}))=
v_n(\bar\phi(F_0(a_0,a_1)),\dots,\bar\phi(F_n(a_0,\dots,a_{n+1})))$
for every $(a_0,\dots,a_{n+1})\in W_{n+1}(A)$.
\end{claim}
\begin{pfclaim} Directly from \eqref{eq_functor_F-V} we derive :
$$
(\bar\phi(F_0(a_0,a_1)),\dots,\bar\phi(F_n(a_0,\dots,a_{n+1})))=
(F_0(\bar\phi(a_0),\bar\phi(a_1)),\dots,F_n(\bar\phi(a_0),
\dots,\bar\phi(a_n))).
$$
Hence, it suffices to show that
$$
\theta\circ v_{n+1}(b_0,\dots,b_{n+1})=
v_n(F_0(b_0,b_1),\dots,F_n(b_0,\dots,b_{n+1}))
$$
for every $(b_0,\dots,b_{n+1})\in W_{n+1}(A/J_0)$. The latter
identity follows from \eqref{eq_basic-F}.
\end{pfclaim}

(iii): We take $u:=\lim_{n\in\N}u_n$, where the maps $u_n$
are as in (ii), for every $n\in\N$. With this definition,
it is clear that $u$ is continuous for the topology
$\cT_{W(R)}$. Next we remark that
$\pi_0\circ u\circ\tau_R=\bar\phi\circ\bomega_0\circ\tau_R=\bar\phi$,
hence :
\set\begin{equation}\label{eq_identity-u}
u\circ\tau_R=\phi
\end{equation}
by lemma \ref{lem_liftings}(iii.a). Now \eqref{eq_explicit}
holds by \eqref{eq_new-form}, \eqref{eq_identity-u} and the
continuity of $u$.
\end{proof}

\begin{remark}\label{rem_tensor-Witt}
(i)\ \
Let $\cN_p$ be the full subcategory of $\Z\tdu\TopAlg$
whose objects are the discrete topological rings in which
$p$ is nilpotent. For every perfect topological $\F_p$-algebra
$E$ and every topological ring $W$, we consider the functors
$$
\bar H_E,H_W:\cN_p\to\Set
$$
such that
$$
\bar H_E(A):=\Hom_{\Z\tdu\TopAlg}(E,A/pA)
\qquad
H_W(A):=\Hom_{\Z\tdu\TopAlg}(W(E),A)
$$
for every $A\in\Ob(\cN_p)$ (where $A/pA$ is endowed
with the discrete topology); to any continuous ring
homomorphism $A\to A'$ of objects of $\cN_p$, the
functor $\bar H_E$ assigns the mapping
$$
\bar H_E(\phi):H_E(A)\to H_E(A')
\qquad
(\bar\psi:E\to A/pA)\mapsto(\phi\otimes_\Z\F_p)\circ\bar\psi
$$
and likewise, $H_W(\phi):H_W(A)\to H_W(A')$ assigns
to any $\psi:A\to A'$ the map $\phi\circ\psi$. Then
proposition \ref{prop_lift-Witt}(iii) says that we have
a natural isomorphism of functors
\set\begin{equation}\label{eq_characterize}
H_{W(E)}\isom\bar H_E
\end{equation}
that assigns to every $A\in\Ob(\cN_p)$ the bijection
$H_{W(E)}(A)\isom\bar H_E(A)$ : $\phi\mapsto\phi\otimes_\Z\F_p$.

(ii)\ \
Suppose that $E$ is a perfect topological $\F_p$-algebra
whose topology is linear, complete and separated. Then
the existence of an isomorphism of functors
\eqref{eq_characterize} characterizes $W(E)$ up to natural
isomorphism, in the category of complete and separated
topological rings whose topology is linear and coarser
than the $p$-adic topology. Indeed, suppose that $W'$
is any other such topological ring, endowed with an
isomorphism of functors $H_{W'}\isom\bar H_E$. There
follows an isomorphism of functors $\omega:H_{W(E)}\isom H_W$.
Now, let $I\subset W(E)$ be any open ideal; the quotient
$W(E)/I$ is an object of $\cN_p$, and the projection
$\pi_I:W(E)\to W(E)/I$ is an element of $H_{W(E)}(W(E)/I)$,
so we get a continuous ring homomorphism
$\omega(\pi_I):W'\to W(E)/I$. Moreover, if $J\subset I$
is another open ideal, the projection
$\pi_{I,J}:W(E)/J\to W(E)/I$ is a morphism of $\cN_p$,
and we have
$$
\pi_{I,J}\circ\omega(\pi_J)=H_{W'}(\pi_{I,J})(\omega(\pi_J))
=\omega(H_{W(E)}(\pi_{I,J})(\pi_J))=\omega(\pi_I)
$$
by the naturality of $\omega$. Hence, we have a well
defined inverse system $(\omega(\pi_I)~|~I\subset W(E))$
indexed by the cofiltered system of open ideals of $W(E)$;
the limit of this system is a well defined continuous
ring homomorphism $\alpha:W'\to W(E)$. Swapping the roles
of $W'$ and $W(E)$, the inverse natural transformation
$\omega^{-1}:H_W\isom H_{W(E)}$ yields likewise a continuous
ring homomorphism $\beta:W(E)\to W'$, as $W'$ is complete
and separated by assumption.

It remains to check that $\alpha$ and $\beta$ are mutually
inverse isomorphisms. To this aim, let $I'\subset W'$ be
any open ideal, $\pi'_{I'}:W'\to W'/I'$ the projection,
and set $I:=\Ker\,(\omega^{-1}(\pi'_{I'}):W(E)\to W'/I')$,
so that $\omega^{-1}(\pi'_{I'})$ factors as the composition
of $\pi_I$ and a continuous ring homomorphism
$\tau:W(E)/I\to W'/I'$. Notice that
$$
\begin{aligned}
\pi'_{I'}\circ\beta\circ\alpha=\,&
\omega^{-1}(\pi'_{I'})\circ\alpha=
\tau\circ\pi_I\circ\alpha=
\tau\circ\omega(\pi_I)=
H_{W'}(\tau)(\omega(\pi_I)) \\
=&\,\omega(H_{W(E)}(\tau)(\pi_I))=\omega(\tau\circ\pi_I)
=\omega(\omega^{-1}(\pi'_{I'}))=\pi'_{I'}
\end{aligned}
$$
which implies that $\beta\circ\alpha=\one_{W'}$. Likewise
one checks that $\beta\circ\alpha=\one_{W(E)}$, whence the
contention.

(iii)\ \
Let $E$, $E'$ and $E''$ be three perfect topological
$\F_p$-algebras whose topologies are linear, complete
and separated, and $E\to E'$, $E\to E''$ two continuous
ring homomorphisms. It is clear that
$(E'\otimes_EE'',\cT^\otimes_{E',E''})$ is a perfect topological
$\F_p$-algebra (notation of \eqref{subsec_tensor-topol}),
hence the same holds for $E'\,\hat\otimes_EE''$ (example
\ref{ex_discrete-Witt}(ii)), and in light of
\eqref{subsec_tensor-topol} and (i), we get natural
isomorphisms of functors
$$
H_{W(E'\,\hat\otimes_EE'')}\isom H_{W(E')}\times_{H_{W(E)}}H_{W(E'')}
\xleftarrow{\sim}H_{W(E')\,\hat\otimes_{W(E)}W(E'')}
$$
whence, by (ii), a natural isomorphism of topological rings
\set\begin{equation}\label{eq_tensor-Witt}
W(E'\,\hat\otimes_EE'')\isom W(E')\,\hat\otimes_{W(E)}W(E'').
\end{equation}
\end{remark}

\begin{proposition}\label{prop_top-basis}
Let $A$ be any discrete $\F_p$-algebra. We have :
\begin{enumerate}
\item
$W(A)$, endowed with its $p$-adic topology, is a flat
and topologically free $\Z_p$-module.
\item
If $A$ is perfect, and $(g_\lambda~|~\lambda\in\Lambda)$
is any basis of the $\F_p$-vector space underlying $A$,
then $(\tau_A(g_\lambda)~|~\lambda\in\Lambda)$ is a
topological basis of the topological $\Z_p$-module
underlying $W(A)$.
\end{enumerate}
\end{proposition}
\begin{proof} The element $p$ is regular in $W(A)$,
by proposition \ref{prop_reduced-Witt}(i), so $W(A)$
is torsion-free, hence flat, as a $\Z_p$-module, for
its natural $\Z_p$-module structure as in example
\ref{ex_discrete-Witt}(iv). In order to check both
(i) and (ii), it then suffices to show the following
general

\begin{claim}\label{cl_top-basis}
Let $R$ be any artinian ring, $\kappa$ (resp. $\fm$) 
the residue field (resp. the maximal ideal) of $R$,
and $M$ any flat $R$-module.
Let also $(f_\lambda~|~\lambda\in\Lambda)$ be any system
of elements of $M$, whose image in $M_0:=M\otimes_R\kappa$
is a basis for the latter $\kappa$-vector space. Then
the induced $R$-linear map
$$
\psi:R^{\oplus\Lambda}\to M
\qquad
e_\lambda\mapsto f_\lambda
\qquad
\text{for every $\lambda\in\Lambda$}
$$
is an isomorphism (here $(e_\lambda~|~\lambda\in\Lambda)$
is the canonical basis of $R^{\oplus\Lambda}$).
\end{claim}
\begin{pfclaim}[] Endow $R^{\oplus\Lambda}$ and $M$
with the $\fm$-adic filtrations, and let
$\gr_\bullet R^{\oplus\Lambda}$ and $\gr_\bullet M$ be
the associated graded $\kappa$-vector spaces. By virtue
of \cite[Th.22.3]{Mat}, the natural $\kappa$-linear map
$$
\fm^n/\fm^{n+1}\otimes_\kappa M_0\to\gr_n M
$$
is an isomorphism for every $n\in\N$. This easily implies
that the $\kappa$-linear map
$\gr_\bullet:\gr_\bullet R^{\oplus\Lambda}\to\gr_\bullet M$
induced by $\psi$, is an isomorphism for every $n\in\N$.
Then the assertion follows from
\cite[Ch.III, \S2, n.8, Cor.3]{BouAC}.
\end{pfclaim}
\end{proof}

\sset\subsubsection{}\label{subsec_combinatorial-Teich}
We shall need a combinatorial identity that expresses
the Teichm\"uller representative of a finite sum of
elements of a perfect topological $\F_p$-algebra $A$,
as a series of terms, each of which is itself a
Teichm\"uller representative.
To state the result, let 
$$
\Sigma^{(k)}_n:=
\{(\sigma_0,\dots,\sigma_k)\in
p^{-n}\N^{k+1}\setminus p^{1-n}\N^{k+1}
~|~\sigma_0+\cdots+\sigma_k=1\}
\qquad
\text{for every $k,n\in\N$}.
$$
and set $\Sigma^{(k)}:=\bigcup_{n\in\N}\Sigma^{(k)}_n$
for every $k\in\N$. Also for every
$a:=(a_0,\dots,a_k)\in A^{k+1}$ and every
$\sigma:=(\sigma_0,\dots,\sigma_k)\in\Sigma^{(k)}$ set
$$
a^\sigma:=a_0^{\sigma_0}\cdots a_k^{\sigma_k}
$$
(notice that the fractional powers $a_i^{\sigma_i}$ are
well defined, since $A$ is perfect).

\begin{proposition}\label{prop_combinatorial}
For every integer $k\in\N$ the following holds :
\begin{enumerate}
\item
There exists a mapping
$\Sigma^{(k)}\to\Z_p\ :\ \sigma\mapsto c_\sigma$ such that
$$
\tau_A(a_0+\cdots+a_k)=\sum_{n\in\N}p^n\cdot
\sum_{\sigma\in\Sigma^{(k)}_n}c_\sigma\cdot\tau_A(a^\sigma)
$$
for every perfect\/ $\F_p$-algebra $A$ and every
$a:=(a_0,\dots,a_k)\in A^{k+1}$.
\item
$c_\sigma=1$ for every $\sigma\in\Sigma_0^{(k)}$.
\end{enumerate}
\end{proposition}
\begin{proof} To begin with, let us show the following :

\begin{claim}\label{cl_tau-expansion}
With the notation of the proposition, for every $n\in\N$ we have :
$$
(\tau_A(a_0)^{1/p^n}+\cdots+\tau_A(a_k)^{1/p^n})^{p^n}\equiv
\tau_A(a_0+\cdots+a_k)
\pmod{p^nW(A)}.
$$
\end{claim}
\begin{pfclaim} Since $\tau_A$ is a splitting of
$\bomega_0$, for every $n\in\N$ we have
$$
\tau_A(a_0)^{1/p^n}+\cdots+\tau_A(a_k)^{1/p^n}\equiv
\tau_A(a_0+\cdots+a_k)^{1/p^n}
\pmod{pW(A)}
$$
and then the claim follows immediately from
lemma \ref{lem_basic-cong}(i).
\end{pfclaim}

Now, consider first the case where
$A=\F_p[T_0^{1/p^\infty},\dots,T_k^{1/p^\infty}]$ is the
universal perfect $\F_p$-algebra in the variables
$T_0,\dots,T_k$. We endow $A$ with the discrete topology.
Set $T_\bullet:=(T_0,\dots,T_k)$, and notice that the
system of monomials $(T_\bullet^\mu~|~\mu\in\N[1/p]^{k+1})$
is a basis for the $\F_p$-vector space underlying $A$, so
the system $(\tau_A(T_\bullet)^\mu~|~\mu\in\N[1/p]^{k+1})$
is a topological basis of the topologically free
$\Z_p$-module underlying $W(A)$ (proposition
\ref{prop_top-basis}(ii)). Hence, the image of
the same system in $W(A)/p^nW(A)$ is a basis for
the latter free $\Z/p^n\Z$-module, for every $n\in\N$.
Consider then the image $S_n$ of $S:=\tau_A(T_0+\cdots+T_k)$
in $W(A)/p^nW(A)$; from claim \ref{cl_tau-expansion}
we see that $S_n$ is uniquely a sum of products
$c^{(n)}_\mu T_\bullet^\mu$ with $c^{(n)}_k\in\Z/p^n\Z$
and $T_\bullet^\mu$ a monomial of the above type, with
the further condition that
$$
c^{(n)}_\mu=0
\qquad\text{unless}\qquad
\mu\in\bigcup_{i=0}^n\Sigma^{(k)}_i.
$$
Moreover, the uniqueness of the resulting expression
means that the image of $c^{(m)}_\mu$ in $\Z/p^n\Z$
agrees with $c^{(n)}_\mu$, for every such $\mu$ and
every $m\geq n$. Summing up, we get the sought identity
for $S$, and (ii) follows as well. Next, if $A$ is any
perfect $\F_p$-algebra, there exists a unique (continuous)
map
$$
\psi:\F_p[T_0^{1/p^\infty},\dots,T_k^{1/p^\infty}]\to A
$$
of $\F_p$-algebras such that $\psi(T_i)=a_i$ for
$i=0,\dots,k$. Since the sought identity is already
known for $T_\bullet$, we deduce it for $a$, by
virtue of \eqref{eq_Teich-functorial}.
\end{proof}

\sset\subsubsection{}\label{subsec_annihilator}
Let $A$ be any perfect $\F_p$-algebra,
$\underline b:=(b-n~|~n\in\N)\in W(A)$ any element, and set
$$
\Ann_A(\underline b):=\{a\in A~|~\tau_A(a)\cdot\underline b=0\}.
$$

\begin{corollary}\label{cor_annihilator}
With the notation of \eqref{subsec_annihilator}, the subset
$\Ann_A(\underline b)$ is a radical ideal of $A$.
\end{corollary}
\begin{proof} Let $a\in\Ann_A(\underline b)$ for some $n\in\N$.
We show first that $a^\mu\in\Ann_A(\underline b)$ for every
$\mu\in\N[1/p]$. Indeed, in view of proposition
\ref{prop_Teich-series}(i), the assumption means that
$a^{p^n}b_n=0$ for every $n\in\N$, and therefore
$a^{p^{n-k}}b_n^{p^{-k}}=0$ for every $k,n\in\N$, since $A$
is perfect. Hence $a^{p^{n-k}}b_n=0$ as well, for every
$n,k\in\N$, which means that
$\tau_A(a^{p^{-k}})\cdot\underline b=0$ for every $k\in\N$
(again by proposition \ref{prop_Teich-series}(i)); the
assertion is a straightforward consequence.

It then remains only to check that $\Ann_A(\underline b)$
is an ideal of $A$. However, if $a,a'\in\Ann_A(\underline b)$
and $x\in A$ is any element, it is clear that
$ax\in\Ann_A(\underline b)$, and proposition
\ref{prop_combinatorial} easily implies that
$a+a'\in\Ann_A(\underline b)$ as well, as required.
\end{proof}

In the same vein, though the Teichm\"uller mapping is not
a ring homomorphism, one can attach to it a continuous map
on Zariski spectra, as explained in the following proposition.

\begin{proposition}\label{prop_like-a-ring-map}
Let $E$ be any perfect $\F_p$-algebra; the following holds :
\begin{enumerate}
\item
For every $\fp\in\Spec\,W(E)$, the subset $\tau_E^{-1}\fp$
is a prime ideal of $E$.
\item
The resulting map
$$
\Spec\,\tau_E:\Spec\,W(E)\to\Spec\,E
\qquad
\fp\mapsto\tau_E^{-1}\fp
$$
is continuous and spectral.
\item
Let $\fp\in\Spec\,W(E)$ be any prime ideal, set
$\fq:=\Spec\,\tau_E(\fp)\in\Spec\,E$, and let
$\pi_\fp:W(E)\to\kappa(\fp)$, $\pi_\fq:E\to\kappa(\fq)$
be the projections. If\/ $\fp$ is a closed subset in the
$p$-adic topology of\/ $W(E)$, then $\tau_E$ induces a
morphism of multiplicative monoids
$$
\tau_\fp:\kappa(\fq)\to\kappa(\fp)
\qquad\text{such that}\qquad
\tau_\fp\circ\pi_\fq=\pi_\fp\circ\tau_E.
$$
\item
Let $E'$ be any other perfect $\F_p$-algebra, and
$\phi:E\to E'$ any continuous ring homomorphism. Then
the resulting diagram commutes :
$$
\xymatrix{ \Spec\,W(E') \ar[rr]^-{\Spec\,\tau_{E'}}
\ar[d]_{\Spec\,W(\phi)} & & \Spec\,E' \ar[d]^{\Spec\,\phi} \\
\Spec\,W(E) \ar[rr]^-{\Spec\,\tau_E} & & \Spec\,E.
}$$
\end{enumerate}
\end{proposition}
\begin{proof}(i): Say that $x,y\in E$ and $\tau_E(xy)\in\fp$;
since $\tau_E$ is a multiplicative map, it follows easily
that either $x$ or $y$ lies in $\tau_E^{-1}\fp$. By the same
token, if $x\in\tau^{-1}_E\fp$ and $y\in E$, then
$xy\in\tau^{-1}_E\fp$.
Next, suppose that $x,y\in\tau^{-1}_E\fp$; by proposition
\ref{prop_combinatorial}, we may write $\tau_E(x+y)$ as
a $p$-adically convergent series $\sum_{n\in\N}p^n\cdot z_n$,
where each $z_n$ is a finite $\Z_p$-linear combination of
terms of the form $\tau_E(x^\lambda y^{1-\lambda})$, with
$\lambda,1-\lambda\in\N[1/p]$. It follows easily that
$z_n=\tau_E(x^{1/p})\cdot z'_n+\tau_E(y^{1/p})\cdot z''_n$,
for suitable $z'_n,z''_n\in W(E)$. Summing up, we get
$$
\tau_E(x+y)=a\cdot\tau_E(x^{1/p})+b\cdot\tau_E(y^{1/p})
\qquad\text{where}\qquad
a:=\sum_{n\in\N}c_nz'_n
\qquad
b:=\sum_{n\in\N}c_nz''_n.
$$
However, by the foregoing we already know that
$x^{1/p},y^{1/p}\in\tau^{-1}_E\fp$, so $\tau_E(x+y)\in\fp$,
whence the contention.

(ii): Let $f\in E$ any element; directly from the
definition we get
\set\begin{equation}\label{eq_was-directly}
(\Spec\,\tau_E)^{-1}(\Spec\,E_f)=\Spec\,A_{\tau_E(f)}
\end{equation}
which says that $\Spec\,\tau_E$ is continuous and spectral.

(iii): We check first that $\tau_E$ descends to a morphism
$\bar\tau:E/\fq\to W(E)/\fp$ of multiplicative monoids. To
this aim, let $x\in E$ and $y\in\fq$ be any two elements;
arguing as in the proof of (i), we see that
$d:=\tau_E(x+y)-\tau_E(x)$ can be written as $p$-adically
convergent series $\sum_{n\in\N}d_n$, where
$d_n\in p^n\tau_E(y^{1/p^n})W(E)$ for every $n\in\N$. We have
already remarked that $\tau_E(y^{1/p^n})\in\fp$ for every
$n\in\N$, and since $\fp$ is $p$-adically closed in $W(E)$,
it follows that $d\in\fp$ as well, whence the claim.
Now, by construction $\bar\tau^{-1}(0)=\{0\}$, so $\bar\tau$
extends uniquely to a morphism of monoids $\tau_\fp$ as sought.

(iv) follows immediately from \eqref{eq_Teich-functorial}.
\end{proof}

\sset\subsubsection{}\label{subsec_mon-fract-powers}
For the further investigation of the ring of Witt vectors
$W(A)$ of a given perfect topological $\F_p$-algebra $A$,
it is useful to consider special types of ideals of $A$
and $W(A)$ defined by certain combinatorial conditions
that we proceed to explain.

$\bullet$\ \
First, consider any {\em $p$-perfect} monoid $P$,
{\em i.e.} such that the $p$-Frobenius endomorphism
$\Phi_P$ of $P$ is bijective. A basic example of
$p$-perfect monoid is $\N[1/p]:=S^{-1}_p\N$, the
localization of $\N$ at its multiplicative subset
$S_p:=\{p^k~|~k\in\N\}$.
Notice that this condition on $P$ means that for every
$x\in P$ and every $\lambda\in\N[1/p]$, the fractional
power $x^\lambda$ is well defined : indeed, we may write
$\lambda=p^{-n}a$ for some $a,n\in\N$ and set
$x^\lambda:=\Phi_P^{-n}(x^a)$. It is easily seen that this
definition is independent of the choice of $a$ and $n$.

$\bullet$\ \
For any ideal $I\subset P$ and every $\lambda=p^{-n}a\in\Z[1/p]$
we define the {\em angular power} :
$$
I^{\La\lambda\Ra}:=\bigcup_{r\in\N}\Phi_P^{-n-r}(I^{p^ra})
\qquad\text{if $\lambda>0$, and}\qquad
I^{\La\lambda\Ra}:=P
\qquad
\text{if $\lambda\leq 0$}.
$$
Notice that $I^{\La\lambda\Ra}$ is an increasing union of
ideals of $P$, and a simple inspection shows that the
definition does not depend on the choice of $a$ and $n$.
Explicitly, if $\lambda\geq 0$, the ideal $I^{\La\lambda\Ra}$
is the subset of all products of the form
$a_1^{\mu_1}\cdots a_k^{\mu_k}$, where $k\in\N$ is any
integer, $a_1,\dots,a_k$ are arbitrary elements of $I$,
and $\mu_1,\dots,\mu_k$ are arbitrary rational numbers
in $\N[1/p]$ such that $\sum_{i=1}^k\mu_i=\lambda$.
Moreover, if $S\subset I$ is any system of generators
for $I$, then $I^{\La\lambda\Ra}$ is generated by the set
of all products as above, where $a_1,\dots,a_k\in S$.

\begin{lemma}\label{lem_mon-fract-powers}
In the situation of \eqref{subsec_mon-fract-powers},
the following holds :
\begin{enumerate}
\item
If $(I_j~|~j\in J)$ is any filtered family of ideals of
$P$, we have
$$
(\bigcup_{j\in J}I_j)^{\La\lambda\Ra}=
\bigcup_{j\in J}I_j^{\La\lambda\Ra}
\qquad
\text{for every $\lambda\in\N[1/p]$}.
$$
\item
For every ideal $I\subset P$, every $n\in\N$ and every
$\lambda,\mu\in\N[1/p]$ we have
\begin{enumerate}
\item
$(I^n)^{\La\lambda\Ra}=(I^{\La\lambda\Ra})^n=I^{\La n\lambda\Ra}$.
\item
$I^{\La\lambda\Ra}I^{\La\mu\Ra}=I^{\La\lambda+\mu\Ra}$.
\item
$(I^{\La\lambda\Ra})^{\La\mu\Ra}=I^{\La\lambda\mu\Ra}$.
\end{enumerate}
\item
$I^{\La\lambda\Ra}J^{\La\lambda\Ra}=(IJ)^{\La\lambda\Ra}$
for every ideals $I,J\subset P$ and every $\lambda\in\N[1/p]$.
\item
Suppose that $I\subset P$ is a non-empty finitely generated
ideal generated by $r$ elements of $P$. Then
$I^{\La\lambda\Ra}\subset I^n$ for every $n\in\N$ such that
$\lambda\geq n+r-1$.
\end{enumerate}
\end{lemma}
\begin{proof}(i): This is clear from the explicit
description of the angular power furnished in
\eqref{subsec_mon-fract-powers}.

(ii.a): Indeed, say that $\lambda=p^ka$ for some $a,n\in\N$;
we have
$$
(I^n)^{\La\lambda\Ra}=\bigcup_{r\in\N}\Phi_P^{-r-k}(I^{nap^r})=
I^{\La n\lambda\Ra}=\bigcup_{r\in\N}
(\Phi_P^{-r-k}(I^{ap^r})^n=(I^{\La\lambda\Ra})^n.
$$

(ii.b): Clearly
$I^{\La\lambda\Ra}I^{\La\mu\Ra}\subset I^{\La\lambda+\mu\Ra}$.
For the converse inclusion, say that
$a_1^{t_1}\cdots a_k^{t_k}\in I^{\La\lambda+\mu\Ra}$, for
given $a_1,\dots,a_k\in I$ and rational numbers
$t_1,\dots,t_k$ in $\N[1/p]$ such that
$\sum_{i=1}^kt_i=\lambda+\mu$. If $\mu=0$, there is
nothing to prove; otherwise, let $j<k$ be the largest
integer such that $s:=\sum_{i=0}^jt_i\leq\lambda$.
Then $t_{j+1}>t'_{j+1}:=\lambda-s$, and if we set
$t''_{j+1}:=t_{j+1}-t'_{j+1}$ we get
$$
a_1^{t_1}\cdots a_{j+1}^{t'_{j+1}}\in I^{\La\lambda\Ra}
\qquad
a_{j+1}^{t''_{j+1}}a_{j+2}^{t_{j+2}}\cdots a_k^{t_k}\in I^{\La\mu\Ra}
$$
whence the contention.

(ii.c): Say that $\lambda=p^{-n}a$ and $\mu=p^{-m}b$;
taking into account (ii.a) we may compute
$$
(I^{\La\lambda\Ra})^{\La\mu\Ra}=
\bigcup_{r\in\N}\Phi_P^{-m-r}((I^{\La\lambda\Ra})^{bp^r})=
\bigcup_{r\in\N}\Phi_P^{-m-r}((I^{bp^r})^{\La\lambda\Ra})=
\bigcup_{r\in\N}\Phi_P^{-m-r}(\bigcup_{s\in\N}\Phi_P^{-n-s}(I^{abp^{r+s}}))
$$
whence the sought identity.

(iii): Indeed, for $\lambda=p^{-n}a$ we may compute
$$
I^{\La\lambda\Ra}J^{\La\lambda\Ra}=
\bigcup_{r\in\N}\Phi_P^{-n-r}(I^{ap^r})\cdot
\bigcup_{r\in\N}\Phi_P^{-n-r}(J^{ap^r})=
\bigcup_{r\in\N}\Phi_P^{-n-r}(I^{ap^r}J^{ap^r})
$$
whence the contention.

(iv): Let $(x_1,\dots,x_r)$ be such a system of
generators for $I$, so that $I^{\La\lambda\Ra}$
is generated by all the products of the form
$x:=x_1^{t_1}\cdots x_r^{t_r}$, where
$t_1,\dots,t_r\in\N[1/p]$ and $t_1+\cdots+t_r=\lambda$.
Then, for every such product and every $i=1,\dots,r$,
let $a_i\in\N$ be the unique integer such that
$t_i\in[a_i,a_i+1[$, and set $s:=a_1+\cdots+a_r$;
clearly $x\in I^s$. Suppose now that $\lambda\geq n+r-1$;
if $t_i\in\N$ for every $i=1,\dots,r$, it follows that
$s=\lambda\geq n$, since $r>0$. Lastly, if there exists
$i\leq r$ such that $t_i\not\in\N$, then $s+r>\lambda$,
whence $s\geq n$ again, and the assertion follows.
\end{proof}

\begin{remark}\label{rem_Witt-are-f-adic}
(i)\ \
Keep the situation of \eqref{subsec_mon-fract-powers},
and suppose moreover that $P$ is a submonoid of the
multiplicative monoid $(A,\cdot)$ of a given ring $A$.
Then, for every ideal $I\subset P$ and every
$\lambda\in\Z[1/p]$ we may consider the ideal
$I^{\La\lambda\Ra}A$ of the ring $A$. In this
situation we also define :
$$
I^{\lfloor r\rfloor}A:=\bigcap_{\substack{ \mu\in\Z[1/p] \\
                          \mu<r\ \ }} I^{\La\mu\Ra}A
\qquad
I^{\lceil s\rceil}A:=\bigcup_{\substack{ \mu\in\Z[1/p] \\
                          \mu>s\ \ }} I^{\La\mu\Ra}A
\qquad
\text{for every $r\in\R$}.
$$
From lemma \ref{lem_mon-fract-powers}(ii.a,ii.c) we deduce
immediately :
$$
(I^{\La\lambda\Ra}A)^n=(I^{\La\lambda\Ra})^nA=I^{\La n\lambda\Ra}A
\qquad
I^{\La\lambda\Ra}I^{\La\mu\Ra}A=I^{\La\lambda+\mu\Ra}A
\qquad
I^{\La\lambda\Ra}J^{\La\lambda\Ra}A=(IJ)^{\La\lambda\Ra}A
$$
for every ideals $I,J\subset P$ and every $\lambda,\mu\in\N[1/p]$
and every $n\in\N$.

(ii)\ \
Especially, if $A$ is a perfect $\F_p$-algebra, we may
take $P:=(A,\cdot)$. In this case, notice that every
ideal $J$ of the ring $A$ is also an ideal of the underlying
monoid $P$, so the angular powers of $J$ are well defined
ideals of $P$; moreover, we easily see that
$$
(IA)^{\La\lambda\Ra}A=I^{\La\lambda\Ra}A
\qquad
\text{for every $\lambda\in\N[1/p]$ and every ideal $I\subset P$}.
$$
Combining with lemma \ref{lem_mon-fract-powers}(ii.c) we get
$$
(I^{\La\lambda\Ra}A)^{\La\mu\Ra}A=I^{\La\lambda\mu\Ra}A
\qquad
\text{for every $\lambda\in\N[1/p]$ and every ideal $I\subset P$}.
$$
Furthermore, suppose that $I$ is an ideal of the ring $A$
generated (in the usual ring-theoretic sense) by a subset
$S\subset A$ consisting of $r$ elements, and let $J_S$ be
the ideal of $P$ generated by $S$; then $I^n=J^n_SA$ and
$I^{\La\lambda\Ra}A=J^{\La\lambda\Ra}_SA$ for every $n\in\N$
and $\lambda\in\N[1/p]$, so lemma \ref{lem_mon-fract-powers}(iv)
yields
$$
I^{\La\lambda\Ra}A\subset I^n
\qquad
\text{for every $n\in\N$ and every $\lambda\in\N[1/p]$
such that $\lambda\geq n+r-1$}.
$$

(iii)\ \
If $A$ is as in (ii), we may also take
$$
P:=\Img\,\tau_A\subset W(A)
$$
which is a $p$-perfect submonoid of the multiplicative monoid
$(W(A),\cdot)$. Clearly $\tau_A$ induces an isomorphism of
monoids $(A,\cdot)\isom P$; thus, if $a_\bullet:=(a_j~|~j\in J)$
is any system of elements of $A$ generating an ideal $I$ of the
monoid $(A,\cdot)$, we may regard $I$ as an ideal of $P$, and we
shall also use the notation
$$
[a_\bullet]^{\La\lambda\Ra}:=I^{\La\lambda\Ra}W(A).
$$

(iv)\ \
For instance, if $A$ is as in (ii), from lemma
\ref{lem_mon-fract-powers}(ii.a) and \eqref{eq_new-form}
we see that
\set\begin{equation}\label{eq_angular-descript}
W(I^{\La\lambda\Ra},r)=
\biggl\{\sum_{n=0}^\infty p^n\cdot\tau_A(b_n)~|~
b_n\in I^{\La\lambda\Ra}A\ \text{for every $n\leq r$}\biggr\}
\end{equation}
for every ideal $I\subset A$, every $\lambda\in\N[1/p]$ and
every $r\in\N$ (notation of remark \ref{rem_Witt-limit}(iv)),
and likewise we can describe $W(I^{\La\lambda\Ra})$. Moreover,
for every such $I$, notice that
\set\begin{equation}\label{eq_two-p-adic-tops}
p^kW(A)\cap W(I^{\La\lambda\Ra})=p^kW(I^{\La\lambda\Ra})
\qquad
\text{for every $\lambda\in\N[1/p]$ and every $k\in\N$}.
\end{equation}
Indeed, say that
$\underline b:=(b_n~|~n\in\N)\in p^kW(A)\cap W(I^{\La\lambda\Ra})$;
in view of lemma \ref{lem_mon-fract-powers}(ii.a), this
means that $b_n=0$ for $n=0,\dots,k-1$, and
$b_n\in I^{\La p^n\lambda\Ra}$ for every integer $n\geq k$.
Clearly
$\underline c:=(b_{n+k}^{1/p^k}~|~n\in\N)\in W(I^{\La\lambda\Ra})$,
and $p^k\cdot\underline c=\underline b$, whence the claim.
As an immediate consequence, we deduce that the $0$-th ghost
map $\bomega_0$ induces an $A$-linear isomorphism
\set\begin{equation}\label{eq_ghost-angular}
W(I^{\La\lambda\Ra})/pW(I^{\La\lambda\Ra})\isom I^{\La\lambda\Ra}A
\qquad
\text{for every $\lambda\in\N[1/p]$}.
\end{equation}
\end{remark}

\begin{lemma}\label{lem_return-to-ic}
Let $A$ be any perfect $\F_p$-algebra, and $I,J\subset A$
two ideals. We have :
\begin{enumerate}
\item
$I^{\lfloor r\rfloor}A=\bar{\mathrm{i.c.}}(I,A,r)$\ \ 
for every $r\in\R_+$ (notation of remark
{\em\ref{rem_add-a-parameter}}).
\item
Suppose that $I$ and $J$ are finitely generated, and
let $(a_1,\dots,a_n)$ (resp. $(b_1,\dots,b_n)$) be a
finite system of generators for $I$ (resp. $J$).
Let also $q>1$ be a rational number with
$$
b_i-a_i\in\mathrm{i.c.}(I,A,q)
\qquad
\text{for every $i=1,\dots,n$}
$$
and suppose furthermore that
\begin{enumerate}
\alphaenu
\item
either, the radical of $I$ equals the radical of $J$
\item
or else, $I$ is contained in the Jacobson radical of $A$.
\end{enumerate}
Then $I^{\La 1\Ra}A=J^{\La 1\Ra}A$.
\end{enumerate}
\end{lemma}
\begin{proof}(i): We are easily reduced to showing that
$$
I^{\La\lambda\Ra}A\subset\mathrm{i.c.}(I,A,\lambda)
\subset I^{\La\lambda'\Ra}A
\qquad
\text{for every $\lambda'\in\N[1/p]$ such that
$\lambda'<\lambda$}.
$$
To this aim, a direct inspection of the definition
shows that
$$
\Phi_A(\mathrm{i.c.}(I,A,\lambda))=\mathrm{i.c.}(I,A,p\lambda)
\qquad\text{and}\qquad
\Phi_A(I^{\La\lambda\Ra}A)=I^{\La p\lambda\Ra}A
\qquad
\text{for every $\lambda\in\N[1/p]$}
$$
so that we may replace $\lambda$ by $p^n\lambda$ for a
suitable $n\in\N$, after which we may assume that
$\lambda\in\N$. Moreover, since the subset
$S:=\{k\lambda/p^n~|~k,n\in\N\}$ is dense in $\N[1/p]$,
we may assume that $\lambda'\in S$. In this case, taking
into account \eqref{eq_lower-m} and lemma
\ref{lem_mon-fract-powers}(ii.a), we may further
replace $I$ by $I^{\La\lambda\Ra}A$, $\lambda$ by $1$, and
$\lambda'$ by $\lambda'/\lambda$, so we reduce to
checking that
$$
I^{\La 1\Ra}A\subset\mathrm{i.c.}(I,A)
\subset I^{\La\mu\Ra}A
\qquad
\text{for every $\mu\in\N[1/p]$ such that $\mu<1$}.
$$
To this aim, say that $x\in I^{\La 1\Ra}$; by
definition, this means that $x^{p^n}\in I^{p^n}$ for
some $n\in\N$, whence $x\in\mathrm{i.c.}(I,A)$.
Next, let $x\in\mathrm{i.c.}(I,A)$, so that
$x^n\in\sum_{i=0}^{n-1}x^iI^{n-i}$ for some $n\in\N$.

\begin{claim}\label{cl_estimate-xI}
$x^m\in\sum_{i=0}^{n-1}x^iI^{m-i}$ for every $m\geq n$.
\end{claim}
\begin{pfclaim} We argue by induction on $m-n$. If
$m=n$, there is nothing to prove, so assume that the
claim is already known for some $m\geq n$. In this
case, we get
$$
x^{m+1}\in\sum_{i=1}^nx^iI^{m+1-i}=x^nI^{m+1-n}+
\sum_{i=1}^{n-1}x^iI^{m+1-i}.
$$
Moreover, $x^nI^{m+1-n}\in\sum^{n-1}_{i=0}x^iI^{m+1-i}$,
so the claim holds for $m+1$.
\end{pfclaim}

It follows from claim \ref{cl_estimate-xI} that
$x^{p^k}\in I^{p^k+1-n}$, whence $x\in I^{\La 1-(n-1)/p^k\Ra}$
for every $k\in\N$ such that $p^k\geq n$, whence the
contention.

(ii): Under the stated conditions, lemma \ref{lem_approx-ic}
says that $\mathrm{i.c.}(I,A)=\mathrm{i.c.}(J,A)$. Together
with \eqref{eq_twice-ic}, we deduce that
$$
\mathrm{i.c.}(I,A,q)=\mathrm{i.c.}(J,A,q).
$$
Hence $J\subset I+\mathrm{i.c.}(I,A,q)\subset I^{\La 1\Ra}A$
and likewise
$I\subset J+\mathrm{i.c.}(J,A,q)\subset J^{\La 1\Ra}A$, whence
the assertion.
\end{proof}

\sset\subsubsection{}\label{subsec_Witt-are-f-adic}
Let $A$ be a perfect $\F_p$-algebra,
$a_\bullet:=(a_1,\dots,a_k)$ a finite system of elements
of $A$, and $I\subset A$ (resp. $\cJ\subset W(A)$) the
ideal generated by $a_\bullet$ (resp. by
$\tau_A(a_1),\dots,\tau_A(a_k)$). Set
$$
\cI:=pW(A)+\cJ
\qquad
W(\lambda):=W(I^{\La\lambda\Ra})
\qquad
\text{for every $\lambda\in\N[1/p]$}
$$
(notation of remark \ref{rem_Witt-limit}(iv)). Lastly,
denote by $\cT_A$ (resp. $\cT_{W(A)}$) the topology of $A$
(resp. of $W(A,\cT_A)$). We complement the generalities
of lemma \ref{lem_mon-fract-powers} with the following:

\begin{proposition}\label{prop_morel}
In the situation of \eqref{subsec_Witt-are-f-adic}, 
the following holds :
\begin{enumerate}
\item
$[a_\bullet]^{\La\lambda\Ra}\subset W(\lambda)\subset
[a_\bullet]^{\La\lambda'\Ra}$
for every $\lambda'<\lambda$ in $\N[1/p]$ (notation
of remark {\em\ref{rem_Witt-are-f-adic}(iii)}).
\item
If $\cT_A$ agrees with the $I$-adic topology, $\cT_{W(A)}$
agrees with the $\cI$-adic topology.
\end{enumerate}
\end{proposition}
\begin{proof}(ii): Notice that $\tau_A(b^\mu)=\tau_A(b)^\mu$
for every $b\in A$ and every $\mu\in\N[1/p]$. Taking into
account \eqref{lem_mon-fract-powers} and proposition
\ref{prop_combinatorial}, we deduce that
$$
W(\lambda,r):=W(I^{\La\lambda\Ra},r)=
p^{r+1}W(A)+[a_\bullet]^{\La\lambda\Ra}
\qquad
\text{every $r\in\N$}
$$
which, in light of lemma \ref{lem_mon-fract-powers}(iv),
shows that the linear topology induced by the cofiltered
system of ideals $(W(\lambda,r)~|~r\in\N,\ \lambda\in\N[1/p])$
agrees with the $\cI$-adic topology. But by the same token
we see that the latter family is a fundamental system of
open neighborhoods of zero in $W(A)$, whence the contention.

(i): Since $W(\lambda,r)$ is an open ideal in the $p$-adic
topology of $W(A)$ for every $r\in\N$, we see that
$W(\lambda)$ is a closed ideal in $W(A)$ for the $p$-adic
topology of $W(A)$. On the other hand, from 
\eqref{eq_two-p-adic-tops} it follows that the $p$-adic
topology of $W(\lambda)$ agrees with the topology
induced from the $p$-adic topology of $W(A)$. Taking into
account proposition \ref{prop_Witt-is-complete}(ii), we
deduce that $W(\lambda)$ is complete and separated for its
$p$-adic topology, and then \eqref{eq_ghost-angular} and
claim \ref{cl_top-basis} imply that $W(\lambda)$ is a
topologically free $\Z_p$-module for this topology, and
any lifting of an $\F_p$-basis of $I^{\La\lambda\Ra}A$
yields a topological basis. However, clearly we may find :
\begin{itemize}
\item
a subset $S\subset\N[1/p]^{\oplus k}$ consisting of sequences
$(\mu_1,\dots,\mu_k)$ such that $\mu_1+\cdots+\mu_k=\lambda$
\item
a system of elements $(b_\mu~|~\mu\in S)$ of $A$ such
that the family $(b_\mu\cdot a^\mu~|~\mu\in S)$
is an $\F_p$-basis of $I^{\La\lambda\Ra}A$, where
$a^\mu:=a_1^{\mu_1}\cdots a_k^{\mu_k}$ for every $\mu\in S$
\end{itemize}
so that the system $(\tau_A(b_\mu\cdot a^\mu)~|~\mu\in S)$
is a topological $\Z_p$-basis of $W(\lambda)$. Now, clearly
$[a_\bullet]^{\La\lambda\Ra}\subset W(\lambda)$; to check the
second inclusion, let $\underline w\in W(\lambda)$ be any
element. By construction, there exists a unique system
$(d_\mu~|~\mu\in S)$ of elements of $\Z_p$ such that
$$
\underline w=\sum_{\mu\in S}d_\mu\cdot\tau_A(b_\mu\cdot a^\mu)
$$
where the series converges in the $p$-adic topology of
$W(A)$, and the set $\{\mu\in S~|~d_\mu\notin p^n\Z_p\}$
is finite for every $n\in\N$. Pick any $N\in\N$ such that
$\lambda-\lambda'\geq kp^{-N}$, and for every $r\in\R$
define
$$
\bar r\in p^{-N}\Z
\qquad
r^*\in[0,p^{-N}[
\qquad\text{fulfilling the condition :}\qquad
r=\bar r+r^*.
$$
Then, for every $\mu\in S$ let also
$\bar\mu:=(\bar\mu_1,\dots,\bar\mu_k)$,
$\mu^*:=(\mu_1^*,\dots,\mu_k^*)$, and set
$$
S':=\{\sigma\in p^{-N}\N^{\oplus k}
      ~|~\lambda\geq\sigma_1+\cdots+\sigma_k>\lambda'\}.
$$
With this notation, it is easily seen that $\bar\mu\in S'$
for every $\mu\in S$. Hence, let also
$$
c_\sigma:=\sum_{\substack{ \mu\in S \\ 
                          \bar\mu=\sigma}}
d_\mu\cdot\tau_A(b_\mu\cdot a^{\mu^*})
\qquad
\text{for every $\sigma\in S'$}
$$
(notice that this series converges $p$-adically in $W(A)$,
so $c_\sigma$ is well defined); we get
$$
\underline w=\sum_{\sigma\in S'}c_\sigma\cdot\tau_A(a^\sigma).
$$
To conclude, it suffices to remark that $S'$ is a finite
set, and $a^\sigma\in I^{\La\lambda'\Ra}A$ for every $\sigma\in S'$.
\end{proof}

\begin{example}\label{ex_Witt-ring}
Let $A:=\F_p[T^{1/p^\infty}]$, endow $A$ with its $T$-adic
topology $\cT_T$, and let $(A^\wedge,\cT^\wedge_T)$ be the
completion of $(A,\cT_T)$. According to lemma
\ref{lem_Witt-limit}(iv), the topological ring
$W(A^\wedge,\cT^\wedge_T)$ is the completion of $W(A,\cT_T)$.
In light of example \ref{ex_discrete-Witt}(v) and
proposition \ref{prop_morel}(ii), the latter is
naturally isomorphic to the $(p,T)$-adic completion of
$\Z_p\{T^{1/p^\infty}\}$, or equivalently, the $(p,T)$-adic
completion of $\Z_p[T^{1/p^\infty}]$. Explicitly, this is
the ring of all power series \eqref{eq_fract-powerseries}
where $(\lambda_n~|~n\in\N)$ is any sequence of elements
of $\N[1/p]$, and $(a_n~|~n\in\N)$ is any sequence
of elements of $\Z_p$, with the property that
$$
\lim_{n\to\infty}v_p(a_n)+\lambda_n=\infty
$$
(where $v_p:\Z_p\to\N\cup\{\infty\}$ is the $p$-adic
valuation).
\end{example}

\sset\subsubsection{}\label{subsec_seminorm-on-W}
Let now $A$ be any topological ring, $|\cdot|_A$ a
real-valued semi-norm on $A$ (see definition
\ref{def_valuation}(v)). For any non-empty subset
$S\subset A$ we let
$$
|S|_A:=\sup_{a\in S}|a|_A\in\R_+\cup\{+\infty\}.
$$
Also, for every real number $\rho\in[0,1]$ we consider
the mapping
$$
|\cdot|_\rho:W(A)\to\R_+\cup\{+\infty\}
\qquad
(a_n~|~n\in\N)\mapsto
\sup_{n\in\N}|a_n|_A^{p^{-n}}\cdot\rho^n.
$$
For instance, notice that $|0|_\rho=0$ and $|1|_\rho=1$;
more generally, we have $|\tau_A(a)|_\rho=|a|_A$ for every
$a\in A$. We let
$$
W(A,\rho):=\{\underline a\in W(A)~|~|\underline a|_\rho\in\R_+\}.
$$

\begin{lemma}\label{lem_semi-norm}
With the notation of \eqref{subsec_seminorm-on-W}, the
following holds :
\begin{enumerate}
\item
$W(A,\rho)$ is an additive subgroup of\/ $W(A)$ and
both $W(A,0)$ and $W(A,1)$ are subrings of\/ $W(A)$.
\item
$|\underline a-\underline b|_\rho\leq
\max(|\underline a|_\rho,|\underline b|_\rho)$
\ \ for every $\underline a,\underline b\in W(A)$.
\item
$|\underline a\cdot\underline b|_{\rho_1\rho_2}\leq
|\underline a|_{\rho_1}\cdot|\underline b|_{\rho_2}$
\ \ for every $\underline a,\underline b\in W(A)$
and every $\rho_1,\rho_2\in[0,1]$.
\item
$\rho^p\cdot|F_A(\underline a)|_{\rho^p}\leq|\underline a|^p_\rho$
\ \ for every $\underline a\in W(A)$ and every $\rho\in[0,1]$.
\item
$|V_A(\underline a)|_{\rho^{1/p}}=
\rho^{1/p}\cdot|\underline a|_\rho^{1/p}$
\ \ for every $\underline a\in W(A)$ and every $\rho\in[0,1]$.
\item
$|\tau_A(a)\cdot\underline b|_\rho\leq|a|_A\cdot|\underline b|_\rho$
\ \ for every $a\in A$ and every $\underline b\in W(A)$, and if\/
$|\cdot|_A$ is a valuation, this inequality is actually
an equality.
\end{enumerate}
\end{lemma}
\begin{proof} Let
$\underline a:=(a_n~|~n\in\N),\underline b:=(b_n~|~n\in\N)$
be any two elements of $W(A)$, and suppose that
$M:=\max(|\underline a|_\rho,|\underline b|_\rho)\in\R_+$.
This condition means that
$$
|a_n|_A,|b_n|_A\leq(M/\rho^n)^{p^n}
\qquad
\text{for every $n\in\N$}.
$$
Set $\underline c:=\underline a+\underline b$,
$\underline d:=\underline a\cdot\underline b$ and
recall that $\underline c=(c_n~|~n\in\N)$ (resp.
$\underline d=(d_n~|~n\in\N)$), where
$c_n=S_n(a_0,\dots,a_n,b_0,\dots,b_n)$
(resp. $d_n=P_n(a_0,\dots,a_n,b_0,\dots,b_n)$ :
notation of \eqref{subsec_Witt-laws}) for every
$n\in\N$.
By construction, $S_n$ is a $\Z$-linear combination
of monomials $Q(X_0,\dots,X_n,Y_0,\dots,Y_n)$ of total
degree $p^n$ (remark \ref{rem_homogeneous-laws}(i));
we claim that
$$
r:=|Q(a_0,\dots,a_n,b_0,\dots,b_n)|_A\leq(M/\rho^n)^{p^n}
\qquad
\text{for every such $Q$}.
$$
Indeed, say that $Q=\prod_{i=0}^nX_i^{n_i}Y_i^{m_i}$; then
$\sum_{i=0}^n(n_i+m_i)p^i=p^n$, and therefore
$$
r\leq\prod_{i=0}^n(M/\rho^i)^{(n_i+m_i)p^i}=M^{p^n}\cdot\rho^{-t}
\qquad
\text{where $t:=\sum_{i=0}^ni\cdot(n_i+m_i)p^i\leq np^n$}.
$$
Since $\rho\leq 1$, the claim follows. Lastly, arguing
with the polynomials $I_n$ of \eqref{subsec_Witt-laws}
that are bihomogeneous of degree $(p^n,0)$, we get
as well $|-\underline a|_\rho\leq|\underline a|_\rho$,
and therefore $|-\underline a|_\rho=|\underline a|_\rho$
for every $\underline a\in W(A)$. Assertion (ii) is an
immediate consequence, and we also deduce that $W(A,\rho)$
is an additive subgroup of $W(A)$.

(iii): Likewise, $P_n$ is a $\Z$-linear combination of
monomials $R$ of bidegree $(p^n,p^n)$, and a similar
calculation shows that
$$
|R(a_0,\dots,a_n,b_0,\dots,b_n)|_A\leq
|\underline a|_{\rho_1}\cdot|\underline b|_{\rho_2}
\cdot(\rho_1\rho_2)^{-np^n}
\qquad
\text{for every such $R$}
$$
(details left to the reader), whence the assertion.
Together with (ii), this also implies that $W(A,0)$
and $W(A,1)$ are subrings of $W(A)$, so the proof of
(i) is complete.

(iv) is proven in the same way ; for every $n\in\N$
one estimates $|F_n(a_0,\dots,a_{n+1})|$ by decomposing
$F_n$ as a sum of bihomogeneous monomials of degree
$(p^{n+1},0)$ : the details shall be left to the reader.

(v): Set $M:=|\underline a|_\rho$, fix a real number
$\eps>0$, and pick $n\in\N$ such that
$|a_n|_A^{p^{-n}}\cdot\rho^n>M-\eps$. Then we have
$$
|a_{k-1}|_A\leq(M/\rho^{k-1})^{p^{k-1}}=
(M^{1/p})^{p^k}\cdot\rho^{-(k-1)p^{k-1}}
\qquad
\text{for every $k>0$}
$$
which is equivalent to
$$
|a_{k-1}|_A^{p^{-k}}\cdot\rho^{k/p}\leq(\rho M)^{1/p}
\qquad
\text{for every $k>0$}.
$$
On the other hand, we have
$|a_n|^{p^{-n-1}}_A\cdot\rho^{(n+1)/p}>(\rho(M-\eps))^{1/p}$
whence the contention.

(v) follows directly from proposition
\ref{prop_Teich-series}(i).
\end{proof}

\begin{proposition}\label{prop_semi-norm-on-W}
In the situation of \eqref{subsec_seminorm-on-W}, suppose
that $A$ is an $\F_p$-algebra. We have:
\begin{enumerate}
\item
$|F_A(\underline a)|_{\rho^p}=|\underline a|_\rho^p$
\ \ for every $\underline a\in W(A)$ and every $\rho\in[0,1]$.
\item
$W(A,\rho)$ is a subring of\/ $W(A)$, and
$|\cdot|_\rho$ restricts to a semi-norm on $W(A,\rho)$.
\item
If moreover $|\cdot|_A$ is a valuation, the same holds
for the restriction of\/ $|\cdot|_\rho$ to $W(A,\rho)$.
\item
If $A$ is perfect and $|\cdot|_A$ is a valuation, we have
$$
\Bigl|\sum_{n\in\N}p^n\cdot\tau_A(b_n)\Bigr|_\rho=
\sup_{n\in\N}|b_n|_A\cdot\rho^n
\qquad
\text{for every sequence $(b_n~|~n\in\N)$ of elements of $A$}. 
$$
\end{enumerate}
\end{proposition}
\begin{proof} Assertion (i) (resp. (iv)) follows easily
from \eqref{eq_simple-F} (resp. from \eqref{eq_new-form}).

(iii): For $\rho=0$, we have $|\underline a|_\rho=|a_0|_A$
for every $\underline a:=(a_n~|~n\in\N)\in W(A)$, so the
assertion is clear in this case. Hence, suppose that
$\rho\in]0,1]$, and let
$\underline a,\underline a'\in W(A,\rho)$ be any two elements.
It is clear that $\underline a\in W(A,\rho-\eps)$ for
every $\eps\in[0,\rho]$, and moreover
\set\begin{equation}\label{eq_semi-limit}
\lim_{\eps\to 0}|\underline a|_{\rho-\eps}=|\underline a|_\rho.
\end{equation}
Furthermore, it is easily seen that, if $\eps>0$, there
exist minimal $N,N'\in\N$ such that
$$
|a_N|^{p^{-N}}_A\cdot(\rho-\eps)^N=|\underline a|_{\rho-\eps}
\qquad\text{and}\qquad
|a'_{N'}|^{p^{-N'}}_A\cdot(\rho-\eps)^{N'}<|\underline a'|_{\rho-\eps}.
$$
We have to check that $|\underline a\cdot\underline a'|_\rho=
|\underline a|_\rho\cdot|\underline a'|_\rho$, and in light
of \eqref{eq_semi-limit} it suffices to show that
$|\underline a\cdot\underline a'|_{\rho-\eps}=
|\underline a|_{\rho-\eps}\cdot|\underline a'|_{\rho-\eps}$
for every $\eps\in[0,\rho]$. We may then replace $\rho$ by
$\rho-\eps$, and assume from start that $N,N'$ are the minimal
integers such that
$$
|a_N|^{p^{-N}}_A\cdot\rho^N=|\underline a|_\rho
\qquad\text{and}\qquad
|a'_{N'}|^{p^{-N'}}_A\cdot\rho^{N'}<|\underline a'|_\rho.
$$
With this notation, write
$\underline a=\underline b+V_A^N(\underline c)$, with
$\underline b:=(a_0,\dots,a_{N-1},0,\dots)$,
$\underline c:=(a_N,a_{N+1},\dots)$. By remark
\ref{rem_semi-norm}(iii) we have
$$
|\underline b|_\rho<|\underline a|_\rho
\qquad\text{and therefore}\qquad
|\underline a|_\rho=|V_A^N(\underline c)|_\rho
$$
from which we deduce that
$|\underline b\cdot\underline a'|_\rho<
|\underline a|_\rho\cdot|\underline a'|_\rho$,
so -- again by remark \ref{rem_semi-norm}(iii) --
we are reduced to checking that
$|V_A^N(\underline c)|_\rho\cdot|\underline a'|_\rho=
|V_A^N(\underline c)\cdot\underline a'|_\rho$.
But proposition \ref{prop_V_A-and-F_A}(ii) and
lemma \ref{lem_semi-norm}(iv) yield
$$
|V_A^N(\underline c)\cdot\underline a'|_\rho=
|V_A^N(\underline c\cdot F^N_A(\underline a'))|_\rho=
|\underline c\cdot F^N_A(\underline a')|^{1/p^N}_{\rho^{p^N}}
\cdot\rho^N
$$
and on the other hand
$$
|V_A^N(\underline c)|_\rho=
|\underline c|^{1/p^N}_{\rho^{p^N}}\cdot\rho^N
\qquad\text{and}\qquad
|F^N_A(\underline a')|^{1/p^N}_{\rho^{p^N}}=
|\underline a'|_\rho
$$
by lemma \ref{lem_semi-norm}(iv) and (i). Say that
$\underline a'':=F^N_A(\underline a')$; then notice
that -- in view of \eqref{eq_simple-F} -- the integer
$N'$ is still the smallest one such that
$$
|a''_{N'}|^{p^{-N'}}_A\cdot\rho^{N'}=|\underline a''|_\rho.
$$
Summing up, we may replace $\underline a$ and $\underline a'$
by $\underline c$ and respectively $\underline a''$,
and $\rho$ by $\rho^N$, and assume from start
that $N=0$, {\em i.e.} that $|a_0|_A=|\underline a|_\rho$.
Likewise, by repeating the same considerations on $\underline a'$,
we may assume that $|a'_0|_A=|\underline a'|_\rho$.
Lastly, notice that
$$
|\underline a\cdot\underline a'|_\rho\geq|a_0a'_0|_A=
|a_0|_A\cdot|a'_0|_A=|\underline a|_\rho\cdot
|\underline a'|_\rho
$$
and the converse inequality is already known, by (ii),
whence the contention.

(ii): Let $\underline a:=(a_n~|~n\in\N)$ and
$\underline b:=(b_n~|~n\in\N)$ be any two elements of
$W(A,\rho)$; we have to show that
$|\underline a\cdot\underline b|_\rho\leq|
\underline a|_\rho\cdot|\underline b|_\rho$. However,
from propositions \ref{prop_Teich-series}(iii) and
\ref{prop_Witt-is-complete}(iv) we get
$$
\underline a\cdot\underline b\equiv
\sum_{i+j=n}V_A^{i+j}(\tau_A(a^{p^j}b^{p^i}))
\pmod{V^{n+1}_A}
\qquad
\text{for every $n\in\N$}
$$
and in view of lemma \ref{lem_semi-norm}(ii) we are then
easily reduced to the case where $\underline a=V^i_A(\tau_A(x))$
and $\underline b=V^j_A(\tau_A(y))$ for some $i,j\in\N$ and
$x,y\in A$. We compute :
$$
\begin{aligned}
|V^{i+j}_A(\tau_A(x^{p^j}y^{p^i}))|_\rho=\, &
\rho^{i+j}\cdot|\tau_A(x^{p^j}y^{p^i})|^{1/p^{i+j}}_{\rho^{p^{i+j}}}
& & \quad\text{(by lemma \ref{lem_semi-norm}(v))} \\
=\, & \rho^{i+j}\cdot|x^{p^j}y^{p^i}|^{1/p^{i+j}}_A  \\
\leq\, & \rho^{i+j}\cdot|x|_A^{1/p^i}\cdot|y|^{1/p^j}_A \\
=\, & |V^i_A(\tau_A(x))|_\rho\cdot|V^j_A(\tau_A(y))|_\rho
& & \quad\text{(again by lemma \ref{lem_semi-norm}(v))}
\end{aligned}
$$
as required.
\end{proof}

\begin{remark} In spite of the inequalities of lemma
\ref{lem_semi-norm}(vi) and proposition
\ref{prop_semi-norm-on-W}(iv), the map $\tau_A:A\to W(A,\rho)$
is not necessarily continuous for the topology on $W(A,\rho)$
induced by the semi-norm $|\cdot|_\rho$, unless $\rho=0$.
\end{remark}

\begin{lemma}\label{lem_new-powers}
Let $E$ be a perfect topological $\F_p$-algebra,
$I\subset E$ any ideal, and denote by $\cR_E$ the
set of all real-valued valuations on $E$ with $|E|=1$
(notation of \eqref{subsec_seminorm-on-W}). For every
$r,r'\in\R_+$ we have :
\begin{enumerate}
\item
$|I^{\lfloor r\rfloor}|_E=|I|_E^r$\ \ for every $|\cdot|_E$
in $\cR_E$.
\item
$(I^{\lfloor r\rfloor})^{\lfloor r'\rfloor}=I^{\lfloor rr'\rfloor}$.
\item
$I^{\lfloor r\rfloor}=(I^{\lfloor r\rfloor})^{\La 1\Ra}$.
\item
$W(I^{\lfloor r\rfloor}E)=\{\underline a\in W(E)~|~
|\underline a|_1\leq|I|^r\ \ \text{for every $|\cdot|$ in $\cR_E$}\}$.
\end{enumerate}
\end{lemma}
\begin{proof}(i): Clearly $|I^{\lfloor r\rfloor}|_E\leq|I|^r_E$.
For the converse inequality, fix any real number
$\eps>0$, and pick $y\in I$ such that $|y|^r_E>|I|^r_E-\eps$.
Since $\N[1/p]$ is dense in $\R_+$, we may find
$\lambda\in\N[1/p]$ such that $\lambda\geq r$ and
$|y|_E^\lambda>|I|^r_E-\eps$. It follows that
$|y^\lambda|=|y|^\lambda\leq|I|^r$ for every $|\cdot|$
in $\cR_E$, and $|y^\lambda|_E>|I|^r_E-\eps$, so that
$|I^{\lfloor r\rfloor}|>|I|^r_E-\eps$, as needed.

(ii): Indeed, by lemmata \ref{lem_strict-closure} and
\ref{lem_return-to-ic}(i), the ideal
$(I^{\lfloor r\rfloor})^{\lfloor r'\rfloor}$ (resp.
$I^{\lfloor rr'\rfloor}$) is the set of all $x\in E$ such that
$|x|\leq|I^{\lfloor r\rfloor}|^{r'}$ (resp. $|x|\leq|I|^{rr'}$)
for every $|\cdot|$ in $\cR_E$.
But from (i) we get $|I^{\lfloor r\rfloor}|^{r'}=|I|^{rr'}$ for
every such $|\cdot|$, whence the assertion.

(iii): This is obvious for $r=0$, and for $r>0$, we have
$I^{\lfloor r\rfloor}\subset(I^{\lfloor r\rfloor})^{\La 1\Ra}
\subset(I^{\lfloor r\rfloor})^{\lfloor 1\rfloor}$
by lemma \ref{lem_return-to-ic}(i), and then the assertion
follows from (ii).

(iv) follows by combining (iii), \eqref{eq_angular-descript}
and proposition \ref{prop_semi-norm-on-W}(iv).
\end{proof}

\subsection{Fontaine rings}
\label{sec_Fontaine-only}
In this section we revisit, and suitably expand, Fontaine's
theory of the {\em fields of norms}, in preparation for its
later employment in the theory of perfectoid rings.

\begin{definition}
(i)\ \
A {\em topological monoid} is the datum $(P,\mu,\cT)$ of a
monoid $(P,\mu)$ and a topology $\cT$ on the set $P$, such
that the composition law $\mu:P\times P\to P$ of $P$ is a
continuous map (where $P\times P$ is endowed with the
product topology). A morphism of topological monoids
$(P,\mu,\cT)\to(P',\mu',\cT')$ is a morphism of monoids
$P\to P'$ that is continuous for the topologies $\cT$ and
$\cT'$. Clearly the topological monoids and their morphisms
form a category
$$
\TopMon.
$$

(ii)\ \
Let $p$ be any prime integer.
We say that $(P,\mu,\cT)$ is a {\em $p$-perfect} topological
monoid if the $p$-Frobenius endomorphism $\bep_P:P\to P$ is
an isomorphism of topological monoids (see definition
\ref{def_exact-phi}(ii)). This means that $P$ is uniquely
$p$-divisible, and $\bep_P$ is a homeomorphism for the
topology $\cT$. We denote by
$$
p\tdu\TopMon
$$
the full subcategory of $\TopMon$ whose objects are
the $p$-perfect topological monoids.
\end{definition}

\begin{remark}\label{rem_TopMon-complete}
Using the criterion of proposition
\ref{prop_complete-criteria}(i) is easily seen that
$\TopMon$ is a complete category : indeed, any family
$((P_i,\cT_i)~|~i\in I)$ of topological monoids admits
a product $P$, which is just the product $P$ of the
underlying topological spaces, endowed with the unique
composition law such that all the projections $P\to P_i$
are morphisms of monoids. Also, if $f,g:(P,\cT)\to(P',\cT')$
is any pair of morphisms in $\TopMon$, the equalizer
of $f$ and $g$ is the equalizer of the underlying morphisms
of monoids, endowed with the topology induced from $\cT$ :
details left to the reader. Likewise, one checks that
$\TopMon$ is also cocomplete.
\end{remark}

\begin{proposition}\label{prop_E-is-right-adj}
The inclusion functor $p\tdu\TopMon\to\TopMon$ admits
both a right and a left adjoint, denoted respectively
$$
\bE:\TopMon\to p\tdu\TopMon
\qquad\text{and}\qquad
\bE^*:\TopMon\to p\tdu\TopMon.
$$
\end{proposition}
\begin{proof} Let $(P,\cT)$ be any topological monoid;
we let $((P_n,\cT_n);\phi_n:P_{n+1}\to P_n~|~n\in\N)$ be
the system of topological monoids such that
$(P_n,\cT_n):=(P,\cT)$ and $\phi_n$ is the $p$-Frobenius
endomorphism of $P$ for every $n\in\N$. In light of remark
\ref{rem_TopMon-complete} we may set
$$
(\bE(P),\cT_{\bE(P)}):=\lim_{n\in\N}\,(P_n,\cT_n)
\qquad
(\bE^*(P),\cT^*_{\bE(P)}):=\colim_{n\in\N}\,(P_n,\cT_n).
$$
So $\bE(P)$ is the set of all sequences $(a_n~|~n\in\N)$
of elements $a_n\in P$ such that $a_n=a_{n+1}^p$ for every
$n\in\N$. Any morphism of topological monoids $f:P\to Q$
induces a morphism of topological monoids
$$
\bE(f):\bE(P)\to\bE(Q)
\qquad
(a_n~|~n\in\N)\mapsto(f(a_n)~|~n\in\N).
$$
Moreover, clearly the Frobenius endomorphism $\bep_{\bE(P)}$
is an automorphism of $\bE(P)$ with inverse $\bff$ given
explicitly by the rule :
$$
\bff(a_n~|~n\in\N)=(a_{n+1}~|~n\in\N)
\qquad
\text{for every $(a_n~|~n\in\N)\in\bE(P)$}.
$$
In other words, $\bff$ is the limit of the system of
morphisms $(\one_P:P_n\to P_{n+1}~|~n\in\N)$, so it
is a continuous map, and therefore $\bep_{\bE(P)}$
is a homeomorphism. Furthermore, the projection to
$P_0$ defines a natural morphism of topological monoids
\set\begin{equation}\label{eq_project-to-P_0}
\bar u_P:\bE(P)\to P
\qquad
(a_n~|~n\in\N)\mapsto a_0
\end{equation}
and if $P$ is already $p$-perfect, then clearly $\bar u_P$
is an isomorphism of topological monoids. Lastly, let
$f:Q\to P$ be any morphism of topological monoids, such that
$Q$ is $p$-perfect. There follows a well defined morphism
$f^\dagger:=\bE(f)\circ\bar u_Q^{-1}:Q\to\bE(P)$ of topological
monoids. Conversely, to any $g:Q\to\bE(P)$ we may attach
the morphism $g^\dagger:=\bar u_P\circ g:Q\to P$.
It is easily seen that the rules $f\mapsto f^\dagger$ and
$g\mapsto g^\dagger$ are mutually inverse bijections, which
concludes the construction of the right adjoint.

Next, let us check that the inverse $\bg$ of $\bep_{\bE^*(P)}$
is continuous. Indeed, any element $\alpha\in\bE^*(P)$ is
the class of a pair $(a,n)$, for some $n\in\N$ and some
$a\in P_n$, and $\bg(\alpha)$ is then the class of
$(a,n+1)$; in other words, $\bg$ is the colimit of the
system of morphisms $(\one_P:P_n\to P_{n+1}~|~n\in\N)$, whence
the contention. This shows that $\bE^*(P)$ is an object
of $p\tdu\TopMon$, and obviously the rule $P\mapsto\bE^*(P)$
yields a well defined functor as sought.
Furthermore, the map $j_0:P\to\bE^*(P)$ defines a natural
transformation which is an isomorphism in case $P$ is
already $p$-perfect. Now, the definition of an adjunction
between $\bE^*$ and the inclusion functor is, {\em mutatis
mutandis}, as in the foregoing case : the details shall
be left to the reader.
\end{proof}

\begin{remark}\label{rem_general}
(i)\ \
By construction, the natural transformation $P\mapsto\bar u_P$
defined in \eqref{eq_project-to-P_0} is the counit of the
adjunction exhibited by proposition \ref{prop_E-is-right-adj}.
Likewise, the unit of this adjunction is the transformation
$Q\mapsto\bar u_Q^{-1}$ for every $p$-perfect topological monoid
$Q$.

(ii)\ \
Taking into account the triangular identities of
\eqref{subsec_adj-pair}, we deduce from (i) the identity
$$
\bE(\bar u_P)=\bar u_{\bE(P)}:\bE(\bE(P))\to\bE(P)
$$
and this map is an isomorphism of $p$-perfect topological
monoids.

(iii)\ \
For every topological monoid $P$, we have
$$
\bE(P)^\times=\bar u_P^{-1}(P^\times).
$$
Indeed, the inclusion
$\bE(P)^\times\subset\bar u_P^{-1}(P^\times)$ is obvious.
For the converse, suppose that $x:=(x_n~|~n\in\N)\in\bE(P)$
and $\bar u_P(x)\in P^\times$. This means that
$x_0\in P^\times$, and since $x_{n+1}^p=x_n$ for
every $n\in\N$, it follows that $x_n\in P^\times$
for every $n\in\N$, so $x\in\bE(P)^\times$.

(iv)\ \
Likewise, it is easily seen that an element $x\in\bE(P)$
is regular if and only if the same holds for $\bar u_P(x)$
(see example \ref{ex_regular}(i)).
\end{remark}

\begin{example} Endow the abelian group $\Q_p/\Z_p$ (resp.
$\Q_p$) with its discrete topology (resp. with its usual
$p$-adic topology); then there exists a natural isomorphism
$\bE(\Q_p/\Z_p)\isom\Q_p$, that identifies the map
$\bar u_{\Q_p/\Z_p}$ with the projection $\Q_p\to\Q_p/\Z_p$
(verification left to the reader).
\end{example}

\sset\subsubsection{}\label{subsec_ring-struct-on-E}
Consider now the category of {\em topological $\F_p$-algebras}
$$
\F_p\tdu\TopAlg
$$
(where $\F_p$ is endowed with its discrete topology : see
definition \ref{def_top-ring}(iii)). If $R$ is an object
of $\F_p\tdu\TopAlg$, the $p$-Frobenius of the
multiplicative monoid $(R,\cdot)$ is the Frobenius
endomorphism $\Phi_R$ of the ring $R$, and a simple
inspection of the proof of proposition
\ref{prop_E-is-right-adj} reveals that the composition
law of $\bE(R,\cdot,\cT)$ is the multiplication map for
a natural structure of topological ring. Moreover, for any
morphism of topological $\F_p$-algebras $f:R\to S$, the
map $\bE(f):\bE(R)\to\bE(S)$ is a morphism of topological
$\F_p$-algebras. In other words, $\bE$ lifts to a functor
\set\begin{equation}\label{eq_E-for-topalg}
\bE:\F_p\tdu\TopAlg\to\F_p\tdu\TopAlg_\mathrm{perf}
\end{equation}
to the full subcategory of $\F_p\tdu\TopAlg$ whose objects
are the {\em perfect} topological $\F_p$-algebras (see
\eqref{subsec_perfect-case}), and \eqref{eq_E-for-topalg}
commutes with the forgetful functors
$(R,+,\cdot,\cT)\mapsto(R,\cdot,\cT)$ to topological
monoids. Furthermore, the map $\bar u_R$ of
\eqref{eq_project-to-P_0} is a morphism of topological
$\F_p$-algebras, so the adjunction exhibited
in the proof of proposition \ref{prop_E-is-right-adj}
lifts to an adjunction between the functor $\bE$ of
\eqref{eq_E-for-topalg} and the forgetful functor
$\F_p\tdu\TopAlg_\mathrm{perf}\to\F_p\tdu\TopAlg$.

The same argument shows that $\bE^*$ lifts to a left
adjoint to the forgetful functor
$$
\bE^*:\F_p\tdu\TopAlg\to\F_p\tdu\TopAlg_\mathrm{perf}.
$$

\begin{remark}\label{rem_topology-of-E}
(i)\ \
Suppose that $(R,\cT)$ is a topological $\F_p$-algebra
whose topology $\cT$ is linear. Then it is easily seen
that a fundamental system of open neighborhoods of
$0\in\bE(R,\cT)$ is given by the intersections of
finitely many elements of the family of ideals
$$
\Phi^n_{\bE(R)}(\bar u{}_R^{-1}(I))
\qquad
\text{for every $n\in\N$ and every open ideal $I\subset R$}.
$$
Especially, the topology $\cT_\bE$ of $\bE(R)$ is also linear.
Moreover, if $\cT$ is separated (resp. and complete)
then the same holds for the topology of $\bE(R)$, since
the latter is a closed subset of $R^\N$, for the product
topology : details left to the reader.

(ii)\ \
In the situation of (i), suppose moreover that $I\subset R$
is an ideal such that $\cT$ agrees with the $I$-adic topology
on $R$, and set $\cI:=\bar u{}_R^{-1}(I)$. Since
$\Phi^n_{\bE(R)}(\cI)\subset\cI^{p^n}\subset\bar u{}_R^{-1}(I^{p^n})$,
we see that $\cT_\bE$ is the linear topology defined by the
system of ideals $(\Phi^n_{\bE(R)}(\cI)~|~n\in\N)$. Especially,
if $\cI$ is finitely generated, $\cT_\bE$ agrees with the
$\cI$-adic topology.

(iii)\ \
Suppose that the Frobenius endomorphism $\Phi_R$ of the
topological $\F_p$-algebra $R$ is surjective. Then
clearly the same holds for $\bar u_R$. Moreover, if
$\Phi_R$ is open and surjective, the same holds for
$\bar u_R$. Indeed, any open subset of $\bE(R)$ is
a union of subsets of the form
$U:=\bE(R)\cap\prod_{n\in\N}U_n$, where each $U_n$ is
open in $R$ and there exists $k\in\N$ such that $U_n=A$
for every $n>k$; with this notation, since $\Phi_R$
is surjective, we have
$$
\bar u_R(U)=
\Phi^k_R(U_k\cap\Phi^{-1}_RU_{k-1}\cap\cdots\cap\Phi^{-k}_RU_0)
$$
and this set is open in $R$, when $\Phi_R$ is an open map.

(iv)\ \
By virtue of proposition \ref{prop_was-get-maddd}(iii),
both the functor $\bE$ for topological monoid and the
one for topological $\F_p$-algebras commute with all
limits. By the same token, the existence of the left
adjoint $\bE^*$ implies that the forgetful functor
$p\tdu\TopMon\to\TopMon$ also commutes with limits,
and likewise for perfect topological $\F_p$-algebras.

(v)\ \
Let $A$ be any $\F_p$-algebra, $I\subset A$ a finitely
generated ideal, and denote by $\cT_I$ (resp. $\cT_d$)
the $I$-adic topology (resp. the discrete topology) on
$A$. Then $(A,\cT_I)$ is a perfect topological $\F_p$-algebra
if and only if the same holds for $(A,\cT_d)$. Indeed,
suppose that the Frobenius endomorphism $\Phi_A$ is
bijective on $A$; in order to show that $\Phi_A$ is
an automorphism of $(A,\cT)$, it suffices to prove that
for every $k\in\N$ there exists $n(k)\in\N$ such that
$I^{n(k)}\subset\Phi(I^n)$; but the latter is clear,
since $I$ is finitely generated (details left to the
reader).
\end{remark}

\begin{theorem}\label{th_fontaine}
Let $(A,\cT_A)$ be any complete and separated topological
ring whose topology $\cT_A$ is linear and coarser than
the $p$-adic topology. Let also $I\subset A$ be a given
ideal that is topologically nilpotent for the topology
$\cT_A$, and $\pi:A\to A/I$ the projection. Endow $A/I$
with the quotient topology $\cT_{A/I}$ induced by $\cT_A$
via $\pi$, and denote by $((A/I)^\wedge,\cT^\wedge_{A/I})$ the
completion of $(A/I,\cT_{A/I})$, and by $j:A\to(A/I)^\wedge$
the completion map. Suppose moreover that :
\begin{enumerate}
\alphaenu
\item
either, $I$ is a closed subset for the topology $\cT_A$
\item
or else, $I$ is finitely generated, and a closed subset
for the $p$-adic topology of $A$.
\end{enumerate}
Then the maps
$$
\bE(\pi):\bE(A,\cT_A)\to\bE(A/I,\cT_{A/I})
\qquad
\bE(j):\bE(A/I,\cT_{A/I})\to\bE((A/I)^\wedge,\cT^\wedge_{A/I})
$$
are isomorphisms of topological monoids.
\end{theorem}
\begin{proof}Let $(J_\lambda~|~\lambda\in\Lambda)$ be a
family of ideals of $A$ which is a fundamental system of
open neighborhoods of $0$ for the topology $\cT_A$.
Suppose first that $I$ is closed for the topology $\cT_A$,
and set $\pi^\wedge:=j\circ\pi:A\to(A/I)^\wedge$; in this
case, $j$ is injective, so the same holds for $\bE(j)$,
and to prove the theorem it suffices to check that
$\bE(\pi^\wedge)$ is an isomorphism. However, we have
natural isomorphisms
$$
A\isom\lim_{\lambda\in\Lambda}\,A/J_\lambda
\qquad
(A/I)^\wedge\isom\lim_{\lambda\in\Lambda}\,A/(I+J_\lambda)
$$
whence, by remark \ref{rem_topology-of-E}(iv), induced
isomorphisms of topological monoids
$$
\bE(A)\isom\lim_{\lambda\in\Lambda}\,\bE(A/J_\lambda)
\qquad
\bE((A/I)^\wedge)\isom
\lim_{\lambda\in\Lambda}\,\bE(A/(I+J_\lambda)).
$$
Thus, it suffices to show that the induced map
$\bE(A/J_\lambda)\to\bE(A/(I+J_\lambda))$ is an
isomorphism of topological monoids for every
$\lambda\in\Lambda$. In this case, we may then
replace $A$ by $A/J_\lambda$, and assume from start
that the topology $\cT_A$ is discrete, hence $I$
is a nilpotent ideal.

\begin{claim}\label{cl_pilpul}
The theorem holds when the topology $\cT_A$
is discrete and $I=pA$.
\end{claim}
\begin{pfclaim} In this case, the $p$-adic filtration
is finite, hence the homomorphism of topological
$\F_p$-algebras $\bar u_{A/pA}:\bE(A/pA)\to A/pA$
lifts to a morphism of monoids $v:\bE(A/pA)\to A$
(lemma \ref{lem_liftings}(iii)). Explicitly, say that
$p^kA=0$ for some $k\in\N$, let
$\bar a_\bullet:=(\bar a_n~|~n\in\N)\in\bE(A/pA)$ be
any element, and pick $a_k\in A$ whose image in $A/pA$
equals $\bar a_k$; then lemma \ref{lem_basic-cong}(i)
shows that $v(\bar a_\bullet)=a_k^{p^k}$. This description
easily implies that $v$ is continuous (for the discrete
topology on $A$). Therefore $v$ corresponds, by adjunction,
to a unique morphism
$v^\dagger:=\bE(v)\circ\bar u^{\ -1}_{\bE(A/pA)}:\bE(A/pA)\to\bE(A)$
of $p$-perfect topological monoids. We have
$$
\bE(\pi)\circ v^\dagger=\bE(\pi\circ v)\circ\bar u^{\ -1}_{\bE(A/pA)}
=\bE(\bar u_{A/pA})\circ\bar u^{\ -1}_{\bE(A/pA)}=\one_{\bE(A/pA)}
$$
by remark \ref{rem_general}(ii). Lastly, the foregoing
explicit description of $v$ implies that
$v^\dagger\circ\bE(\pi)=\one_{\bE(A)}$ (details left to
the reader), and the claim follows.
\end{pfclaim}

Now, we have a commutative diagram of topological monoids
$$
\xymatrix{ \bE(A) \ar[rr]^-{\bE(\pi')} \ar[d]_{\bE(\pi)} & &
\bE(A/pA) \ar[d]^{\bE(\pi'')} \\
\bE(A/I) \ar[rr]^-{\bE(\pi''')} & & \bE(A/(pA+I))
}$$
where $\pi',\pi'',\pi'''$ are the projections. By
claim \ref{cl_pilpul} we know already that $\bE(\pi')$
and $\bE(\pi''')$ are isomorphisms, hence it suffices
to show that the same holds for $\bE(\pi'')$. Thus,
we may further replace $A$ by $A/pA$, and assume
from start that $A$ is also an $\F_p$-algebra. In
this case, the $I$-adic filtration on $A$ fulfills
the conditions of lemma \ref{lem_basic-cong}, hence
lemma \ref{lem_liftings}(iii) applies and says that
the map $\bar u_{A/I}:\bE(A/I)\to A/I$ lifts to a
morphism of monoids $w:\bE(A/I)\to A$, and arguing
as in the proof of claim \ref{cl_pilpul} we see that
$w$ is continuous, so it yields, again by adjunction,
a morphism of topological monoids
$w^\dagger:\bE(A/I)\to\bE(A)$, which is seen to be
an inverse for $\bE(\pi)$, by the same argument as
in the proof of claim \ref{cl_pilpul}. This concludes
the proof, in case $I$ is closed for the topology $\cT_A$.

For the general case, denote by $\bar I\subset A$ the
topological closure of $A$ (for the topology $\cT_A$),
endow $A/\bar I$ with the quotient topology $\cT_{A/\bar I}$
induced by the projection $A\to A/\bar I$, and let
$((A/\bar I)^\wedge,\cT^\wedge_{A/\bar I})$ be the
completion of $(A/\bar I,\cT_{A/\bar I})$; denote also
by $\pi':(A/I)^\wedge\to(A/\bar I)^\wedge$ and
$\bar\pi:A\to(A/\bar I)^\wedge$ the induced maps; we
deduce morphisms of topological monoids
$\bE(A)\to\bE((A/I)^\wedge)\to\bE((A/\bar I)^\wedge)$
whose composition equals $\bE(\bar\pi)$. Notice that
$\bar I$ is still topological nilpotent, hence $\bE(\bar\pi)$
is an isomorphism, by the foregoing case. Also
$\pi'$ is an isomorphism, so the same holds for
$\bE(\pi')$ as well. We set
$v:=\bE(\bar\pi)^{-1}\circ\bE(\pi'):\bE((A/I)^\wedge)\to\bE(A)$.
Clearly $v\circ\bE(\pi)=\one_{\bE(A)}$, so it remains
only to check that $v$ is bijective, or equivalently,
that the same holds for $\bE(\pi)$. To this aim, let
$\cT'_A$ be the $(I+pA)$-adic topology on $A$, and
notice that $\cT'_A$ is still complete and separated
(lemma \ref{lem_fontaine}); moreover, clearly both
$I$ and $pA$ are topologically nilpotent ideals in
$(A,\cT'_A)$. Furthermore, notice that
$$
\bigcap_{n\in\N}(I+(I+pA)^n)=\bigcap_{n\in\N}(I+p^nA)=I
$$
since $I$ is closed in the $p$-adic topology of $A$;
thus, $I$ is closed in the topology $\cT'_A$, so the
foregoing case shows that the map
$\bE(\pi):\bE(A,\cT'_A)\to\bE(A/I)$ is indeed bijective.
\end{proof}

\begin{remark}\label{rem_fontaine}
Keep the situation of theorem \ref{th_fontaine}.

(i)\ \
For every $(\bar a_n~|~n\in\N)\in\bE(A/I)$ and every
$n\in\N$ choose a representative $a_n\in A$ for the class
$\bar a_n$. By inspecting the proof, we see that the inverse
of $\bE(\pi)$ is given by the rule :
\set\begin{equation}\label{eq_u-and-tau}
(\bar a_n~|~n\in\N)\mapsto
(\lim_{k\to\infty} a_{n+k}^{p^k}~|~n\in\N)
\end{equation}
where the convergence is relative to the topology $\cT_A$.

(ii)\ \
Suppose moreover, that $p\in I$; then $\bE(A/I)$ is
a topological $\F_p$-algebra, and the isomorphism
$\bE(\pi)$ can be used to transfer to $\bE(A)$ the
ring structure of $\bE(A/I)$. In light of (i), we get
the following expressions for the addition and
multiplication laws on $\bE(A)$ such that $\bE(\pi)$
becomes an isomomorphism of topological rings :
$$
\begin{aligned}
(a_n~|~n\in\N)+(b_n~|~n\in\N):=\, &
(\lim_{k\to+\infty}(a_{n+k}+b_{n+k})^{p^k}~|~n\in\N) \\
(a_n~|~n\in\N)\cdot(b_n~|~n\in\N):=\, &
(a_nb_n~|~n\in\N)
\end{aligned}
$$
where the limit is computed in the topology $\cT_A$.

(iii)\ \
In the situation of (ii), suppose furthermore that $A$ is
a domain (resp. a field). Then theorem \ref{th_fontaine}
implies immediately that the same holds for $\bE(A/I)$ (and
of course also for $\bE(A)$, for its ring structure as in (ii)).
\end{remark}

\begin{corollary}\label{cor_E-and-completion}
Let $(A,\cT)$ be a perfect topological $\F_p$-algebra,
$I\subset A$ an open and topologically nilpotent ideal,
and denote by $(A^\wedge,\cT^\wedge)$ the separated
completion of $(A,\cT)$. Then there is a natural
isomorphism of topological $\F_p$-algebras
$$
\omega:(A^\wedge,\cT^\wedge)\isom\bE(A/I)
\qquad\text{such that}\qquad
\pi=\omega\circ\bar u_{A/I}
$$
where $\pi:A^\wedge\to A/I$ is the natural map.
\end{corollary}
\begin{proof} Let $I^\wedge$ be the topological closure of
the image of $I$ in $A^\wedge$, and endow $B:=A^\wedge/I^\wedge$
with the quotient topology induced by the projection
$\pi:A^\wedge\to B$; notice that the natural map $A/I\to B$
induces an isomorphism $\phi:A/I\isom B^\wedge$ on the
respective separated completions. Then example
\ref{ex_discrete-Witt}(iii) and theorem
\ref{th_fontaine} imply that the composition
$$
A^\wedge\xrightarrow{\ \bar u^{\ -1}_{A^\wedge}\ }\bE(A^\wedge)
\xrightarrow{\ \bE(\pi)\ }\bE(B)\xrightarrow{\ \bE(\phi)\ }
\bE(A/I)
$$
is an isomorphism, whence the corollary.
\end{proof}

\sset\subsubsection{}\label{subsec_E-of-val.ring}
Let now $(R,|\cdot|_R)$ be any valuation ring with
value group $\Gamma_{\!\!R}$ and residue characteristic $p$.
Notice that the $p$-Frobenius map of $\Gamma_{\!\!R}$ is
injective, so we may regard $\bE(\Gamma_{\!\!R})$ as a
subgroup of $\Gamma_{\!\!R}$. To ease notation, set
$\bE:=\bE(R/pR)$; also, let $\Delta:=\Gamma_{\!\!R}$
in case $R$ is an $\F_p$-algebra, and otherwise let
$\Delta:=O(|p|_R)$ (notation of remark
\ref{rem_ordered-gps}(vii)). We define a mapping
$$
|\cdot|_{\bE}:\bE\to\Gamma_{\!\!\bE\circ}
\qquad
\text{where $\Gamma_{\!\!\bE}:=\bE(\Gamma_{\!\!R})\cap\Delta$}
$$
as follows. Let $\underline a:=(a_n~|~n\in\N)\in\bE$ be any
element, and for every $n\in\N$ choose a representative
$\tilde a_n\in R$ for the class $a_n\in R/pR$. If
$|\tilde a_n|_R\leq|p|_R$ for every $n\in\N$, then we set
$|\underline a|_\bE:=0$; if there exists $n\in\N$
such that $|\tilde a_n|>|p|_R$, then we set
$|\underline a|_\bE:=|\tilde a_n|^{p^n}_R$.
One verifies easily that the definition is independent of
all the choices : the details shall be left to the reader.

\begin{lemma}\label{lem_E-is-val.ring}
In the situation of \eqref{subsec_E-of-val.ring}, 
the pair $(\bE,|\cdot|_\bE)$ is a valuation ring.
\end{lemma}
\begin{proof} The bijection of remark
\ref{rem_valuations}(vii) attaches to $\Delta$ the
prime ideal $\fp(\Delta)=\bigcap_{n\in\N}p^nR$, and
we set $R':=R/\fp(\Delta)$. Since the set of ideals
of $R'$ is naturally identified with the set of those
ideals of $R$ that contain $\fp(\Delta)$, it follows
easily that $R'$ is a valuation ring, and then it is
clear that its valuation is induced by the primary
specialization $|\cdot|_R^\Delta$ of $|\cdot|_R$ (notation
of \eqref{subsec_special-investig}). Notice moreover that
$R'/pR'=R/pR$, so that $\bE(R'/pR')=\bE$.
Summing up, we may replace $R$ by $R'$ after which we
may assume from start that $\Gamma_{\!\!R}=\Delta$, and
that the $p$-adic topology on $R$ is separated. Moreover,
notice that in this case the $p$-adic topology on $R$
agrees with either the valuation topology (see definition
\ref{def_valuation-toplog}) or with the discrete topology;
hence the $p$-adic completion $R^\wedge$ of $R$ is still a
valuation ring (proposition \ref{prop_stays-valuation}(iii)).
Then, since $R^\wedge/pR^\wedge=R/pR$, we may replace
$R$ by $R^\wedge$, and assume from start that $R$ is
$p$-adically complete, in which case remark
\ref{rem_fontaine}(iii) already shows that $\bE$ is
a domain. Moreover, in this situation the isomorphism
$\bE(R)\isom\bE$ of theorem \ref{th_fontaine}
identifies the map $|\cdot|_\bE$ with the mapping
$$
|\cdot|_{\bE(R)}:\bE(R)\to\Gamma_{\!\!R\circ}
\qquad
(a_n~|~n\in\N)\mapsto|a_0|_R.
$$
With this notation, clearly
$|x_\bullet\cdot y_\bullet|_{\bE(R)}=
|x_\bullet|_{\bE(R)}\cdot|y_\bullet|_{\bE(R)}$ for every
$x_\bullet,y_\bullet\in\bE(R)$. By the same token, it
is clear that $x_\bullet$ divides $y_\bullet$ in $\bE(R)$
if and only if $|x_\bullet|_{\bE(R)}\geq|y_\bullet|_{\bE(R)}$ 
so $(\bE,|\cdot|_\bE)$ is a valuation ring.
\end{proof}

The last topics of this section are the graded versions
of the functors $W$ and $\bE$.

\sset\subsubsection{}\label{subsec_graded-Witt}
Let $(\Gamma,+,0)$ be a monoid, $p\in\N$ a prime integer,
$(A,\underline B)$ a topological ring with $\Gamma$-graded
structure, and suppose that the $p$-Frobenius endomorphism
of $\Gamma$ is injective. The assumption implies that the
natural map $\Gamma\to\bE^*(\Gamma)$ is injective (notation
of proposition \ref{prop_E-is-right-adj}), hence we may
regard $\underline B$ naturally as an $\bE^*(\Gamma)$-graded
$\Z$-algebra. We define
$$
W(\underline B):=
\bigoplus_{\gamma\in\bE^*(\Gamma)}\gr_\gamma W(\underline B)
\qquad\text{where}\qquad
\gr_\gamma W(\underline B):=\prod_{n\in\N}\gr_{p^n\gamma}B
\qquad
\text{for every $\gamma\in\bE^*(\Gamma)$}.
$$
From remark \ref{rem_homogeneous-laws}(i) it follows
easily that, for every $\gamma,\mu\in\Gamma$, the addition
and multiplication laws of $W(B)$ restrict respectively to maps
$$
\gr_\gamma W(\underline B)\times\gr_\gamma W(\underline B)\to
\gr_\gamma W(\underline B)
\qquad
\gr_\gamma W(\underline B)\times\gr_\mu W(\underline B)\to
\gr_{\gamma+\mu} W(\underline B).
$$
We may therefore define an addition law on $W(\underline B)$
by taking the direct sum of the addition laws on each
$\Gamma$-graded summand. Likewise, the map $x\mapsto -x$
on $W(B)$ restricts to a map
$\gr_\gamma W(\underline B)\to\gr_\gamma W(\underline B)$
for every $\gamma\in\Gamma$, so that $W(\underline B)$
is an abelian group with this addition law. Next, clearly
$1\in\gr_0W(\underline B)$, so we obtain a ring structure
on $W(\underline B)$, by extending linearly the
multiplication map already defined on each graded term.
Thus, the multiplication map of $W(\underline B)$ agrees
with that of $W(B)$ on each pair of homogeneous elements,
but it will be usually different, on non-homogeneous
elements.

\begin{lemma}\label{lem_graded-Witt}
With the notation of \eqref{subsec_graded-Witt},
the induced map
$$
\phi:W(\underline B)\to W(B)
\qquad
(b_\gamma~|~\gamma\in\bE^*(\Gamma))\mapsto
\sum_{\gamma\in\bE^*(\Gamma)}b_\gamma
$$
is injective (again, notice that $\phi$ agrees with
the inclusion map on homogeneous terms, but it usually
differs from it on non-homogeneous elements).
\end{lemma}
\begin{proof} For every $r\in\N$ and
$\gamma\in\bE^*(\Gamma)$ define
$\gr_\gamma W_r(\underline B)$ as the image of
$\gr_\gamma W(\underline B)$ into $W_r(B)$, and set
$$
W_r(\underline B):=
\bigoplus_{\gamma\in\bE^*(\Gamma)}\gr_\gamma W_r(\underline B).
$$
We shall show, more precisely, that the resulting map
$\phi_r:W_r(\underline B)\to W_r(B)$ is a bijection.
We argue by induction on $r$. For $r=1$, the assertion
follows by a simple inspection. Next, suppose that
$r\geq 0$ and that $\phi_r$ is bijective; we
deduce a commutative ladder with exact rows :
$$
\xymatrix{
0 \ar[r] &
\bigoplus_{\gamma\in\bE^*(\Gamma)}\gr_{p^r\gamma}B
\ar[r] \ar[d]_{\psi_r} &
W_{r+1}(\underline B) \ar[r] \ar[d]_{\phi_{r+1}} &
W_r(\underline B) \ar[r] \ar[d]^{\phi_r} & 0 \\
0 \ar[r] & W_{r+1}(B)\cap\bar V_r \ar[r] &
W_{r+1}(B) \ar[r] & W_r(B) \ar[r] & 0
}$$
where $\bar V_r$ denotes the image of $V_r(B)$
in $W_{r+1}(B)$. It then suffices to check that $\psi_r$
is a bijection; but \eqref{eq_addition-on-Vn} implies
that $\psi_r$ is none else than the restriction of the
additive map $V^r$, and since $\bE^*(\Gamma)$ is
$p$-perfect, the source of $\psi_r$ is just $B$,
whence the contention.
\end{proof}

Due to lemma \ref{lem_graded-Witt}, we may regard
$W(\underline B)$ as a subring of $W(B)$. We also
see that $W_r(B)$ inherits from $\underline B$ a
natural $\bE^*(\Gamma)$-graded $\Z$-algebra structure,
for every $r\in\N$. Moreover,
$\gr_0W(\underline B)=W(\gr_0B)$, so $W(\underline B)$ is
an $\bE^*(\Gamma)$-graded $W(\gr_0B)$-algebra. Set
$$
W(A,\underline B):=(W(A),W(\underline B))
$$
and let also $\bE^*(\Gamma)_{(\bp)}$ be the object of the
category $\bE^*(\Gamma)/\Mnd$ given by the $p$-Frobenius
automorphism
$\bp_{\bE^*(\Gamma)}:\bE^*(\Gamma)\to\bE^*(\Gamma)_{(\bp)}$.

\begin{proposition}
With the notation of \eqref{subsec_graded-Witt}, the
following holds :
\begin{enumerate}
\item
The pair $W(A,\underline B)$ is an $\bE^*(\Gamma)$-graded
structure on $W(A)$ (see definition {\em\ref{def_graded-top}(i)}).
\item
The map $\bomega_0$ is a morphism of topological rings with
$\bE^*(\Gamma)$-graded structures
$$
\bomega_0:W(A,\underline B)\to(A,\underline B).
$$
\item
The Frobenius and Verschiebung maps of $W(A)$ (see
\eqref{subsec_V-and-F}) induce respectively a morphism
of topological rings with\/ $\bE^*(\Gamma)$-graded structures
$$
F_A:W(A,\underline B)\to
(W(A),\bE^*(\Gamma)\times_{\bE^*(\Gamma)_{(\bp)}}W(\underline B))
$$
and a continuous map of\/ $\bE^*(\Gamma)$-graded topological
abelian groups
$$
V_A:\bE^*(\Gamma)\times_{\bE^*(\Gamma)_{(\bp)}}
W(\underline B)\to W(\underline B)
$$
(notation of definition {\em\ref{def_Gamma-graded-algs}(iv)}).
\end{enumerate}
\end{proposition}
\begin{proof}(i): Indeed, a direct inspection shows that
$\gr_\gamma W(\underline B)$ is a closed subset in
$W(A,\cT_A)$, for every $\gamma\in\bE^*(\Gamma)$.
Next, notice that the topology of $W(B)$ agrees with
the topology induced by $W(A)$ via the inclusion map
$W(B)\to W(A)$. Hence, in order to check that the topology
induced by $W(A)$ on $W(\underline B)$ is linear and
defined by a system of open graded ideals, it suffices to
show that the same holds for the topology induced by
$W(B)$. However, let $(J_\lambda~|~\lambda\in\Lambda)$
be any fundamental system of open graded ideals of $B$, and
set $\cJ_{\lambda,r}:=\Ker\,(W(B)\to W_r(B/J_\lambda))$ for
every $\lambda\in\Lambda$ and $r\in\N$; the proof of lemma
\ref{lem_Witt-limit}(iii) shows that
$(\cJ_{\lambda,r}~|~\lambda\in\Lambda;\ r\in\N)$ is a
fundamental system of open ideals of $W(B)$, and we are
reduced to checking that $\cJ_{\lambda,r}\cap W(\underline B)$
is a graded ideal of $W(\underline B)$, for every
$\lambda\in\Lambda$ and $r\in\N$. To this aim, we consider
the commutative diagram
$$
\xymatrix{ W(\underline B) \ar[r]^-\phi \ar[d] & W(B) \ar[d] \\
W_r(\underline B/J_\lambda) \ar[r]^-{\phi_r} & W_r(B/J_\lambda)
}$$
where $\phi_r$ is the isomorphism defined as in the proof
of lemma \ref{lem_graded-Witt}, and the vertical arrows
are the natural projections. We deduce that
$\cJ_{\lambda,r}\cap W(\underline B)$ is the kernel of the
projection $W(\underline B)\to W_r(\underline B/J_\lambda)$,
and the assertion follows easily. It remains to check that
$W(\underline B)$ is a dense subset of $W(A)$, but this is
clear, since $\phi_r$ is bijective for every $r\in\N$.

(ii) follows by a direct inspection of the construction.

(iii) follows easily from remark \ref{rem_homogeneous-laws}(i).
\end{proof}

We have the following graded version of proposition
\ref{prop_lift-Witt} :

\begin{proposition}\label{prop_graded-lift}
Let $p\in\N$ be a prime integer, $\Gamma$ a monoid whose
$p$-Frobenius endomorphism is injective, $(R,\underline S)$
and $(A,\underline B)$ two topological rings with
$\Gamma$-structures, and suppose that:
\begin{enumerate}
\alphaenu
\item
$R$ is a perfect topological $\F_p$-algebra, and
$A$ is complete and separated.
\item
$\underline B$ admits a fundamental system $(J_n~|~n\in\N)$
of graded open ideals such that
$pJ_n+J_n^p\subset J_{n+1}\subset J_n$ for every $n\in\N$.
\end{enumerate}
Then, for every morphism
$\bar\phi:(R,\underline S)\to(A,\underline B)$ of
topological rings with $\Gamma$-graded structures
there exists a unique morphism
$u:W(R,\underline S)\to(A,\underline B)$ of topological
rings with $\Gamma$-graded structures that makes commute
the diagram
$$
\xymatrix{
W(R,\underline S) \ar[rr]^-u \ar[d]_{\bomega_0} & &
(A,\underline B) \ar[d]^\pi \\
(R,\underline S) \ar[rr]^-{\bar\phi} & & (A,\underline B)/J_0
}$$
where $\pi$ is the natural projection.
\end{proposition}
\begin{proof} We have to show that the unique continuous
ring homomorphism $u:W(R)\to A$ provided by proposition
\ref{prop_lift-Witt}(iii) restricts to a map of $\Gamma$-graded
rings $W(\underline S)\to\underline B$. However, recall
that $\gr_\gamma B$ is a closed subset of $A$, so it is
complete and separated for the topology induced by
$\underline B$, for every $\gamma\in\Gamma$, and hence
also for the $p$-adic topology (lemma \ref{lem_fontaine}).
Then, let $\phi:R\to A$ be the unique continuous mapping
characterized as in lemma \ref{lem_liftings}(iii);
by inspecting \eqref{eq_explicit}, we are reduced to
checking that $\phi(\gr_\gamma S)\subset\gr_\gamma B$ for
every $\gamma\in\Gamma$. To this aim, pick any map
$\psi:R\to A$ such that :
\begin{enumerate}
\alphaenu
\item
$\pi\circ\psi=\bar\phi$
\item
$\psi(\gr_\gamma S)\subset\gr_\gamma B$ for every
$\gamma\in\Gamma$
\end{enumerate}
and recall that $\phi(r)=\lim_{n\to+\infty}\psi(r^{p^{-n}})^{p^n}$
for every $r\in R$. The assertion follows easily.
\end{proof}

\sset\subsubsection{}\label{subsec_reintroduce-B*}
Let $\Gamma$ be a monoid, $P$ a $\Gamma$-graded monoid.
Clearly $\bE(P)$ is naturally an $\bE(\Gamma)$-graded monoid.
Moreover, let $\underline B:=(B,\gr_\bullet B)$ be a
$\Gamma$-graded topological $\Z$-algebra; we define
the $\Gamma$-graded monoid
$$
B^*:=\coprod_{\gamma\in\Gamma}\gr_\gamma B
$$
whose composition law is induced by the multiplication map
of $B$. The natural morphism of $\Gamma$-graded monoids
$B^*\to B$ (which is not injective, unless $\Gamma=\{0\}$)
induces a morphism of $\bE(\Gamma)$-graded monoids
$$
\bE(B^*)\to\bE(B).
$$
Furthermore, if $B$ is an $\F_p$-algebra we may define the
$\bE(\Gamma)$-graded subring of $\bE(B)$
$$
\bE(\underline B):=
\bigoplus_{\gamma\in\bE(\Gamma)}\gr_\gamma\bE(B^*)
$$
which we endow with the topology induced from $\bE(B)$, via
the inclusion map $\bE(\underline B)\to\bE(B)$, and the rule
$\underline B\mapsto\bE(\underline B)$ defines a functor from
the category of $\Gamma$-graded topological $\F_p$-algebras
to the category of perfect $\bE(\Gamma)$-graded topological
$\F_p$-algebras. It is easily seen that this functor is right
adjoint to the forgetful functor
$\underline C\mapsto\underline C_{/\Gamma}$ induced by
$\bar u_\Gamma:\bE(\Gamma)\to\Gamma$.

\begin{proposition}\label{prop_gamma-graded-E-funct}
Let $(\Gamma,+,0)$ be a monoid, $(A,\underline B)$ a
topological $\F_p$-algebra with $\Gamma$-graded structure,
such that $A$ is complete and separated, and the $p$-Frobenius
endomorphism $\bp_\Gamma$ of\/ $\Gamma$ is injective. We have :
\begin{enumerate}
\item
$\bE(A,\underline B):=(\bE(A),\bE(\underline B))$ is an
$\bE(\Gamma)$-graded structure on $\bE(A)$. 
\item
The pair $(\bar u_A,\bar u_\Gamma)$ is a morphism of
topological $\F_p$-algebras with graded structures :
$$
\bar u_A:\bE(A,\underline B)\to(A,\underline B).
$$
\item
Let $\Gamma'$ be another monoid whose $p$-Frobenius
endomorphism is injective, and $\Gamma\to\Gamma'$ any
morphism of monoids. Then the identity map of\/ $\bE(A)$
induces an isomorphism of topological rings with
$\Gamma'$-graded structures
$$
\bE(A,\underline B,\Gamma)_{/\Gamma'}\isom
\bE((A,\underline B,\Gamma)_{/\Gamma'}).
$$
\end{enumerate}
\end{proposition}
\begin{proof} Indeed, the definition immediately shows
that $\gr_\gamma\bE(B)$ is closed in the topology of
$\bE(A)$ for every $\gamma\in\Gamma$.
Next, notice that -- under the current assumptions -- we
may regard $\bE(\Gamma)$ as a submonoid of $\Gamma$, and
for every graded ideal $J\subset B$ and every $r\in\N$, set
$$
J^*:=\coprod_{\gamma\in\Gamma}\gr_\gamma J
\qquad\text{and}\qquad
\bE(\underline J):=\bigoplus_{\gamma\in\Gamma}\gr_\gamma\bE(J^*).
$$
So, $\bE(\underline J)$ is a graded ideal of
$\bE(\underline B)$, and from remark
\ref{rem_topology-of-E}(i) it follows that the
topology of $\bE(A)$ induces on $\bE(\underline B)$
the linear topology defined by the cofiltered system of
graded ideals
$(\Phi^r_{\bE(\underline B)}\bE(J)~|~J\subset B;\ r\in\N)$,
where $J$ ranges over all the open graded ideals of $B$.
It remains to check that $\bE(\underline B)$ is dense
in $\bE(A)$. To this aim, let
$\underline a:=(a_n~|~n\in\N)\in\bE(A)$ be any element,
and for every finite subset $S\subset\bE(\Gamma)$
and every $n\in\N$, set $S(n):=\bp^{-n}_\Gamma(S)$; for
every such $S$, we define the element
$$
\underline a_S:=(a_{n,S(n)}~|~n\in\N)\in\bE(\underline B)
$$
where the component
$a_{n,S(n)}\in\bigoplus_{\gamma\in S(n)}\gr_\gamma B$ is
defined as in remark \ref{rem_graded-top-algs}(iii).

\begin{claim}\label{cl_new-claim}
(i)\ \ 
$a_{n,\gamma}=0$ for every $\gamma\in\Gamma\setminus\bE(\Gamma)$
and every $n\in\N$.
\begin{enumerate}
\addenu
\item
$(a_{n,S(n)}~|~S\subset\bE(\Gamma))$ is a Cauchy
net in $B$ whose limit is $a_n$, for every $n\in\N$.
\item
$a_{n,S(n)}=a^p_{n+1,S(n+1)}$, for every $n\in\N$ and every
finite subset $S\subset\bE(\Gamma)$.
\end{enumerate}
\end{claim}
\begin{pfclaim} For every $n\in\N$, let
$(a_{n,\gamma}~|~\gamma\in\Gamma)\in\prod_{\gamma\in\Gamma}B_\gamma$
be the sequence of canonical projections of $a_n$ (see
remark \ref{rem_graded-top-algs}(iii)); in view of lemma
\ref{lem__p-power-Cauchy}(i) we see that $a_{n,p\gamma}=a^p_{n+1,\gamma}$
for every $\gamma\in\Gamma$, and $a_{n,\gamma}=0$ if
$\gamma\in\Gamma\setminus\bp_\Gamma(\Gamma)$. Assertions (i)
and (iii) are immediate consequences.

(ii) follows immediately from (i) and proposition
\ref{prop_Cauchy}(ii).
\end{pfclaim}

From claim \ref{cl_new-claim}(ii,iii) we conclude that the
system $(\underline a_S~|~S\subset\bE(\Gamma))$ is a Cauchy
net in $\bE(\underline B)$ whose limit is $\underline a$, as
required.

(ii) is clear by inspecting the constructions.

(iii): Indeed, by example \ref{ex_monoid-alg-graded}(iii),
any such morphism will be necessarily an isomorphism, and
conversely, consider the natural morphism
$(A,\underline B,\Gamma)\to(A,\underline B,\Gamma)_{/\Gamma'}$
which induces a morphism $\bE(A,\underline B,\Gamma)\to
\bE((A,\underline B,\Gamma)_{/\Gamma'})$ which in turns
factors uniquely through a morphism as sought,
again by example \ref{ex_monoid-alg-graded}(iii).
\end{proof}

\sset\subsubsection{}\label{subsec_E-for-not-F_p-alg}
Let $(A,\underline B)$ be a complete and separated
topological ring with $\Gamma$-graded structure, whose
topology is linear and coarser than the $p$-adic topology.
Endow $A/pA$ with the quotient topology induced by the
projection $\pi_A:A\to A/pA$, and let $(A/pA)^\wedge$ be the
separated completion of $A/pA$; by theorem \ref{th_fontaine},
the map $\pi_A$ induces an isomorphism of topological monoids
$$
\bE(\pi_A):\bE(A)\isom\bE(A/pA)
$$
and the completion map $j:A/pA\to(A/pA)^\wedge$
induces an isomorphism of topological rings
$$
\bE(j):\bE(A/pA)\to\bE((A/pA)^\wedge).
$$
Especially, $\bE(A/pA)$ is a complete and separated
topological ring whose topology is linear (remark
\ref{rem_topology-of-E}(i)), and we may transfer the
ring structure of $\bE(A/pA)$ onto $\bE(A)$ (see remark
\ref{rem_fontaine}(ii)), after which the latter can
also be regarded as a complete and separated topological
ring with a linear topology. Furthermore, suppose that
the $p$-Frobenius map of $\Gamma$ is injective, so that
the topological rings with $\bE(\Gamma)$-structures
$$
((A/pA)^\wedge,\underline B_0):=((A,\underline B)/pB)^\wedge
\qquad\text{and}\qquad
\bE((A/pA)^\wedge,\underline B_0)
$$
are well defined, by proposition
\ref{prop_gamma-graded-E-funct}(i) (and remark
\ref{rem_graded-top-algs}(iv)), and we may transfer
as well this graded structure onto $\bE(A)$ and
$\bE(A/pA)$. We denote by
$$
\bE(A,\underline B)
\qquad\text{and}\qquad
\bE(A/pA,\underline B/p\underline B)
$$
the resulting topological rings with $\bE(\Gamma)$-graded
structures, whose underlying topological rings are respectively
$\bE(A)$ and $\bE(A/pA)$. Define $B'$ as in remark
\ref{rem_graded-top-algs}(ii), and notice that
$pB=pB'\cap B$ and $pA\subset A\cap pB'$, so that
\set\begin{equation}\label{eq_same-p-filtr}
pB=pA\cap B
\end{equation}
and therefore the natural map $\bE(B/pB)\to\bE(A/pA)$
is injective.

\begin{proposition}
In the situation of \eqref{subsec_E-for-not-F_p-alg},
let also $\underline C$ (resp. $\underline C_0$) be the
$\Gamma$-graded dense subring of $\bE(A)$ (resp. of\/
$\bE(A/pA)$). With the notation of \eqref{subsec_reintroduce-B*}
we have :
$$
\underline C_0=\bE(\underline B/p\underline B)
\qquad\text{and}\qquad
\underline C=\bE(\underline B).
$$
\end{proposition}
\begin{proof} Notice that -- as explained in remark
\ref{ex_was-rem-graded-v}(ii) -- since $A$ is complete
and separated, $\gr_\gamma B_0$ is the maximal separated
quotient of $\gr_\gamma(B/pB)$, for every $\gamma\in\Gamma$;
it is also a closed subset of the complete ring
$(A/pA)^\wedge$, so it is the separated completion of
$\gr_\gamma(B/pB)$.

Now, from \eqref{eq_same-p-filtr} we deduce that
the canonical $\gamma$-projection of $(A,\underline B)$
induces a continuous $B_0/pB_0$-linear map
$\pi_\gamma:A/pA\to\gr_\gamma(B/pB)$ for every
$\gamma\in\Gamma$, and in view of the injectivity of
$\bp_\Gamma$, proposition \ref{prop_Cauchy}(ii) and
lemma \ref{lem__p-power-Cauchy}(i) easily imply that
$\pi_{p\gamma}(a^p)=\pi_\gamma(a)^p$ in the
$\Gamma$-graded $\Z$-algebra $B/pB$. To ease notation,
set $D_\gamma:=\gr_\gamma\bE(\underline B)$ and
$D_{0,\gamma}:=\gr_\gamma\bE(\underline B/p\underline B)$
for every $\gamma\in\bE(\Gamma)$; there follows, for every
such $\gamma$, a continuous $\bE(\gr_0B)$-linear map
$$
\pi_{\bE,\gamma}:\bE(A/pA)\to D_{0,\gamma}
\qquad
(a_n~|~n\in\N)\mapsto(\pi_{\gamma/p^n}(a_n)~|~n\in\N).
$$
Let also $\pi^\wedge_\gamma:(A/pA)^\wedge\to\gr_\gamma B_0$
be the completion of $\pi_\gamma$, for every $\gamma\in\Gamma$;
we obtain likewise a continuous $\bE(\gr_0B_0)$-linear map
for every $\gamma\in\bE(\Gamma)$
$$
\pi^\wedge_{\bE,\gamma}:
\bE((A/pA)^\wedge)\to\gr_\gamma\bE(\underline B_0)
\qquad
(a_n~|~n\in\N)\mapsto(\pi^\wedge_{\gamma/p^n}(a_n)~|~n\in\N).
$$
However, $\pi^\wedge_{\bE,\gamma}$ is none else than the
canonical $\gamma$-projection of $\bE((A/pA)^\wedge,\underline B_0)$;
indeed, it is easily seen that it coincides with the latter
map on the dense subring $\bE(\underline B_0)$, and
$\gr_\gamma\bE(\underline B_0)$ is separated, whence the
claim. Moreover, by construction we have a commutative
diagram
$$
\xymatrix{
\bE(A/pA) \ar[rrr]^-{\prod_{\gamma\in\bE(\Gamma)}\pi_{\bE,\gamma}}
\ar[d]_{\bE(j)} & & &
\prod_{\gamma\in\bE(\Gamma)}D_{0,\gamma} \ar[d]^{\prod_{\gamma\in\bE(\Gamma)}i_\gamma} \\
\bE((A/pA)^\wedge) \ar[rrr]^-{\prod_{\gamma\in\bE(\Gamma)}\pi^\wedge_{\bE,\gamma}}
& & &
\prod_{\gamma\in\bE(\Gamma)}\gr_\gamma\bE(\underline B_0)
}$$
where $i_\gamma:D_{0,\gamma}\to\gr_\gamma\bE(\underline B_0)$
is the restriction of $\bE(j)$, for every $\gamma\in\Gamma$.
Then, a simple inspection of this diagram shows that
$\gr_\gamma C_0\subset D_{0,\gamma}$ for every $\gamma\in\Gamma$,
and the converse inclusion is obvious. Lastly, it is
clear that $\bE(\pi_A)(D_\gamma)\subset D_{0,\gamma}$ for every
$\gamma\in\bE(\Gamma)$. To show the converse inclusion, let
$\bar b_\bullet:=(\bar b_n~|~n\in\N)$ be any element
of $D_{0,\gamma}$, and for every $n\in\N$ pick
$b_n\in\gr_{\gamma/p^n}B$ whose image in $\gr_{\gamma/p^n}(B/pB)$
agrees with $\bar b_n$; according to remark
\ref{rem_fontaine}(i), the element
$\bE(\pi_A)^{-1}(\bar b_\bullet)$ is the sequence
$(b'_n~|~n\in\N)$ with $b'_n:=\lim_{k\to\infty}b^{p^k}_{n+k}$
for every $n\in\N$, where the convergence is relative to
the topology of $A$. But since $\gr_\gamma B$ is closed
in $A$ for every $\gamma\in\Gamma$, we see that
$b'_n\in\gr_{\gamma/p^n}B$ for every $n\in\N$, whence
the claim.
\end{proof}

\subsection{Divided power modules and algebras}
This section is a review of the theory of divided
power modules and algebras.

\begin{definition}\label{def_multipolynomial}
Let $A$ be a ring, $k\in\N$ an integer, and
$\underline M:=(M_1,\dots,M_k)$ an object of the abelian
category $(A\Mod)^k$, {\em i.e.} any sequence of $A$-modules.
Let also $Q$ be another $A$-module, and
$\underline d:=(d_1,\dots,d_k)\in\N^{\oplus k}$ any sequence
of $k$ integers.
\begin{enumerate}
\item
We denote by $F_{\underline M},G_Q:A\Alg\to\Set$
the functors given respectively by the rules :
$$
B\mapsto(M_1\times\cdots\times M_k)\otimes_AB
\quad\text{and}\quad
B\mapsto Q\otimes_AB
\qquad
\text{for every $A$-algebra $B$}.
$$
Also, let $M^\vee_i:=\Hom_A(M_i,A)$ for every $i=1,\dots,k$,
and set :
$$
\Sym^\bullet_A(\underline M^\vee):=
\bigotimes_{i=1}^k\Sym^\bullet_A(M_i^\vee)
$$
which we view as a $\N^{\oplus k}$-graded $A$-algebra, with
the grading given by the rule :
$$
\Sym^{\underline n}_A(\underline M^\vee):=
\bigotimes_{i=1}^k\Sym^{n_i}_A(M_i^\vee)
\qquad
\text{for every $\underline n:=(n_1,\dots,n_k)\in\N^{\oplus k}$}.
$$
\item
A {\em homogeneous multipolynomial law} of degree
$\underline d$ on the sequence of $A$-modules $\underline M$,
with values in $Q$, is a natural transformation
$\lambda:F_{\underline M}\to G_Q$, which we write also
$$
\lambda:\underline M\squig Q
$$
such that the following holds. For every $A$-algebra $B$,
every $(x_1,\dots,x_k)\in F_{\underline M}(B)$ and every
$(b_1,\dots,b_k)\in B^k$, we have
$$
\lambda_B(b_1x_1,\dots,b_kx_k)=b_1^{d_1}\cdots b_k^{d_k}\cdot
\lambda_B(x_1,\dots,x_k).
$$
\item
The set of all homogeneous multipolynomial laws
$\underline M\squig Q $ of degree $\underline d$ is denoted
$$
\Pol^{\underline d}_A(\underline M,Q).
$$
\end{enumerate}
\end{definition}

\begin{remark}\label{rem_apriori}
(i)\ \
Let $\underline f:=(f_i:M'_i\to M_i~|~i=1,\dots,k)$
be a sequence of $A$-linear maps, and set
$\underline M':=(M'_1,\dots,M'_k)$; let also $g:Q\to Q'$ be
another $A$-linear map. Suppose that
$\lambda:\underline M\squig Q$ is a homogeneous
multipolynomial law of degree $d$; then we obtain a law
$$
g\circ\lambda\circ\underline f:\underline M'\squig Q'
$$
of the same type, by the rule :
$$
B\mapsto(g\otimes_AB)\circ\lambda_B\circ
((f_1\times\cdots\times f_k)\otimes_AB):
(M'_1\times\cdots\times M'_k)\otimes_AB
\to Q'\otimes_AB
$$
for every $A$-algebra $B$. This shows that the rule
$(\underline M,Q)\mapsto\Pol^{\underline d}_A(\underline M,Q)$
is functorial in both arguments. However, it is not {\em a priori}
obvious that this functor takes values in (essentially) small
sets. The latter assertion will nevertheless be proven shortly.

(ii)\ \
Suppose that $\underline M':=(M'_1,\dots,M'_{k'})$ is
another sequence of $A$-modules, $Q'$ another $A$-module, and
$\underline d':=(d'_1,\dots,d'_{k'})\in\N^{\oplus k'}$ another
sequence of integers. If we have two homogeneous polynomial
laws $\lambda:\underline M\squig Q$ and
$\lambda':\underline M'\squig Q'$ of degrees respectively
$\underline d$ and $\underline  d'$, we can form the tensor
product
$$
\lambda\otimes\lambda':(\underline M,\underline M')\squig
Q\otimes_AQ'
$$
which is the homogeneous multipolynomial law of degree
$(\underline d,\underline d')$ given by the rule :
$$
(\lambda\otimes\lambda')_B(x_1,\dots,x_k,x'_1,\dots,x'_{k'}):=
\lambda_B(x_1,\dots,x_k)\otimes\lambda'_B(x'_1,\dots,x'_{k'})
$$
for every $A$-algebra $B$ and every
$(x_1,\dots,x'_{k'})\in(M_1\times\cdots\times M'_{k'})\otimes_AB$.

(iii)\ \
Let $B$ be any $A$-algebra; if $\lambda:\underline M\squig Q$
is a law as in (i), then the restriction of $\lambda$ to
$B$-algebras yields a law
$\lambda_{/B}:\underline M\otimes_AB\to Q\otimes_AB$. In this
way, we obtain a natural mapping
\set\begin{equation}\label{eq_base-change-pol}
\Pol^{\underline d}_A(\underline M,Q)\to
\Pol^{\underline d}_B(\underline M\otimes_AB,Q\otimes_AB).
\end{equation}
\end{remark}

\begin{lemma}\label{lem_multipolynomial}
In the situation of definition {\em\ref{def_multipolynomial}},
suppose that the $A$-module $M_i$ is free of finite rank, for
every $i=1,\dots,k$. Then, there exists a natural bijection
$$
\Pol^{\underline d}_A(\underline M,Q)\isom
Q\otimes_A\Sym^{\underline d}_A(\underline M^\vee).
$$
\end{lemma}
\begin{proof} Indeed, fix bases $(e_{i1},\dots,e_{ir_i})$ for
each $A$-module $M_i$, and denote by $(e^*_{i1},\dots,e^*_{ir_i})$
the dual basis of $M^\vee_i$. For every $i=1,\dots,k$
there is a natural (injective) map of $A$-algebras
$$
\Sym^\bullet_A(M^\vee_i)\to\Sym^\bullet_A(\underline M^\vee)
\quad :\quad
s\mapsto 1\otimes\cdots\otimes s\otimes\cdots\otimes 1
$$
which allows to view naturally the $e^*_{ij}$ as elements
of $\Sym^\bullet_A(\underline M^\vee)$. Suppose that
$\lambda:\underline M\squig Q$ is a given homogeneous
multipolynomial law of degree $\underline d$, and set
$$
P_\lambda(e^*_{ij}):=\lambda_{\Sym^\bullet_A(\underline M^\vee)}
\Bigl(\sum_{j=1}^{r_1}e_{1j}\otimes e^*_{1j},\dots,
\sum_{j=1}^{r_k}e_{kj}\otimes e^*_{kj}\Bigr)
\in Q\otimes_A\Sym^\bullet_A(\underline M^\vee).
$$
Notice that the terms in parenthesis are elements of
$M_i\otimes_AM^\vee_i\subset M_i\otimes_A\Sym^\bullet_A(M^\vee_i)$
that do not depend on the chosen bases, since they correspond
to the identity automorphisms of $M_i$, under the natural
identification $M_i\otimes_AM^\vee_i\isom\End_A(M_i)$. Now,
let $B$ be any $A$-algebra, and
$(b_{ij}~|~i=1,\dots,k;\ j=1,\dots,r_i)$ any sequence of
elements of $B$. By considering the unique map of $A$-algebras
$\Sym^\bullet_A(\underline M)\to B$ given by the rule : 
$e^*_{ij}\mapsto b_{ij}$ for every $i=1,\dots,k$ and every
$j=1,\dots,r_i$, it is easily seen that
$$
\lambda_B
\Bigl(\sum_{j=1}^{r_1}b_{1j}e_{1j},\dots,\sum_{j=1}^{r_k}b_{kj}e_{kj}\Bigr)
=P_\lambda(b_{ij})\in Q\otimes_AB.
$$
In other words, $\lambda$ is completely determined by
the polynomial $P_\lambda(e^*_{ij})$ with coefficients in $Q$.

Next, set $B:=\Sym^\bullet_A(\underline M^\vee)[Y_1,\dots,Y_k]$;
the homogeneity condition on $\lambda$ implies that
$$
P_\lambda(Y_ie^*_{ij})=\lambda_B
\Bigl(Y_1\cdot\sum_{j=1}^{r_1}e_{1j}\otimes e^*_{1j},\dots,
Y_k\cdot\sum_{j=1}^{r_k}e_{kj}\otimes e^*_{kj}\Bigr)=
Y_1^{d_1}\cdots Y^{d_k}_k\cdot P(e^*_{ij})
$$
which means that
$P_\lambda\in Q\otimes_A\Sym^{\underline d}_A(\underline M^\vee)$.
Conversely, it is clear that any element of
$Q\otimes_A\Sym^{\underline d}_A(\underline M^\vee)$ yields a
homogeneous law of degree $\underline d$, whence the lemma.
\end{proof}

\begin{remark}\label{rem_funct-on-M}
(i)\ \
Lemma \ref{lem_multipolynomial} shows especially that,
if $M_1,\dots,M_k$ are free $A$-modules of finite rank,
the rule $Q\mapsto\Pol_A^{\underline d}(\underline M,Q)$
yields a functor with values in essentially small sets,
and it is clear that the sum of two homogeneous multipolynomial
laws of degree $\underline d$ is still a law of the same type;
likewise, if we multiply such a law by an element of $A$,
we get a new law of the same type. Up to isomorphism, we may
therefore assume that this is a functor $A\Mod\to A\Mod$.

(ii)\ \
Next, set
$$
\Gamma^{\underline d}_A(\underline M):=
(\Sym^{\underline d}_A\underline M^\vee)^\vee.
$$
The lemma can be rephrased by saying that, under the
stated assumptions, the $A$-module
$\Gamma^{\underline d}_A(\underline M)$ represents the
functor
$$
A\Mod\to A\Mod
\quad :\quad
Q\mapsto\Pol_A^{\underline d}(\underline M,Q).
$$

(iii)\ \
The rule $(A,M)\mapsto\Gamma^{\underline d}_A(\underline M)$
is functorial for sequences $\underline M$ of free
$A$-modules of finite rank. Indeed, let
$\underline\phi:\underline M\to\underline N$ be a morphism
of such sequences. For every degree $\underline d$ of length
$k$, the induced map
$\Pol^{\underline d}_A(\underline N,Q)\to
\Pol^{\underline d}_A(\underline M,Q)$ corresponds to a
map
\set\begin{equation}\label{eq_work-it-out}
\Sym^{\underline d}_A(\underline N^\vee)\otimes_AQ\to
\Sym_A^{\underline d}(\underline M^\vee)\otimes_AQ
\end{equation}
that can be worked out as follows. Pick bases
$(e_{ij}~|~j=1,\dots,r_i)$ for $N_i$ and
$(e_{ij}~|~j=1,\dots,s_i)$ for $M_i$, for every $i=1,\dots,k$,
and denote as usual by $(e^*_{ij}~|~j=1,\dots,r_i)$,
$(f^*_{ij}~|~j=1,\dots,s_i)$ the respective dual bases.
Say that $\lambda:N\squig Q$ is a given homogeneous law
of degree $\underline d$, and let $P_\lambda(e^*_{ij})\in
\Sym^{\underline d}_A(\underline N^\vee)\otimes_AQ$
be the corresponding homogeneous polynomial. By
construction, the polynomial $P_{\lambda\circ\underline\phi}
\in\Sym_A^{\underline d}(\underline M^\vee)\otimes_AQ$
corresponding to $\lambda\circ\phi$ is calculated as
$$
\begin{aligned}
P_{\lambda\circ\underline\phi}:=\, &
(\lambda\circ\underline\phi)_{\Sym^\bullet_A(\underline M^\vee)}
\Bigl(\sum_{j=1}^{s_i}f_{ij}\otimes f^*_{ij}~|~i=1,\dots,k\Bigr) \\
=\, & \lambda_{\Sym^\bullet_A(\underline M^\vee)}
\Bigl(\sum_{j=1}^{s_i}\phi_i(f_{ij})\otimes f^*_{ij}~|~i=1,\dots,k\Bigr) \\
=\, & \lambda_{\Sym^\bullet_A(\underline M^\vee)}
\Bigl(\sum_{j=1}^{s_i}e_{ij}\otimes\phi^\vee_i(e^*_{ij})~|~i=1,\dots,k\Bigr) \\
=\, & P_\lambda(\phi^\vee_i(e^*_{ij})).
\end{aligned}
$$
In other words, \eqref{eq_work-it-out} equals
$\Sym^{\underline d}_A(\underline\phi^\vee)\otimes_A\one_Q$,
so if we let
$$
\Gamma^{\underline d}_A(\underline\phi):=
\Sym^{\underline d}_A(\underline\phi^\vee)^\vee
$$
we obtain :
\begin{itemize}
\item
a functor $\Gamma^{\underline d}_A(-)$ from the
full subcategory $(A\Mod^o)^k_\mathrm{ff}$ of $(A\Mod^o)^k$
consisting of all sequences of free $A$-modules of finite
rank, into the category of $A$-modules
\item
for every $A$-module $Q$, a well defined isomorphism of
functors on this subcategory
$$
(Q,\underline M)\mapsto(\Pol^{\underline d}_A(\underline M,Q)
\isom\Hom_A(\Gamma^{\underline d}_A\underline M,Q)
\quad :\quad
A\Mod\times(A\Mod^o)_\mathrm{ff}^k\to A\Mod
$$
\end{itemize}
The following result extends the above observations to
arbitrary sequences $\underline M$ of $A$-modules.
\end{remark}

\begin{theorem}\label{th_construct-Gammas}
Let $\underline M$ be any sequence as in definition
{\em\ref{def_multipolynomial}(i)}. We have :
\begin{enumerate}
\item
$\Pol^{\underline d}_A(\underline M,Q)$ is an
essentially small set, for every $A$-module $Q$.
\item
There exists a functor
\set\begin{equation}\label{eq_Gamma-multi}
\Gamma^{\underline d}_A:(A\Mod)^k\to A\Mod
\end{equation}
with a natural isomorphism of $A$-module valued functors
$$
\Pol^{\underline d}_A(-,Q)\isom
\Hom_A(\Gamma^{\underline d}_A(-),Q)
\qquad
\text{for every $A$-module $Q$}.
$$
\end{enumerate}
\end{theorem}
\begin{proof} Suppose first that all the $A$-modules
$M_1,\dots,M_k$ are finitely presented. Then, for every
$i=1,\dots,k$, pick a surjective $A$-linear map
$f_i:L_i\to M_i$, with $L_i$ free of finite rank.
Next, let $g_i:L_i\oplus L_i\to M_i$ be the map given by
the rule : $(l,l')\mapsto f_i(l-l')$ for every
$(l,l')\in L_i\oplus L_i$; since $M_i$ is finitely presented,
$\Ker\,g_i$ is finitely generated, for every $i=1,\dots,k$,
so we may find a surjective $A$-linear map $L'_i\to\Ker\,g_i$.
By composing with the inclusion $\Ker\,g_i\to L_i\oplus L_i$
and the two projections $L_i\oplus L_i\to L_i$, we deduce a
presentation of $M_i$ :
\set\begin{equation}\label{eq_presentation-M}
\xymatrix{
\Coker\,(h^1_i,h^2_i:L'_i \ar@<.5ex>[r] \ar@<-.5ex>[r]& L_i)\isom M_i
}\end{equation}
such that the induced map
$$
L'_i\to L_i\times_{M_i}L_i
$$
is surjective. Notice that the latter condition is stable
under arbitrary base changes $A\to B$. Set
$\underline L:=(L_1,\dots,L_k)$,
$\underline L':=(L'_1,\dots,L'_k)$ and
$\underline h^j:=(h^j_1,\dots,h^j_k)$ for $j=1,2$.

\begin{claim}\label{cl_presentation-M}
The presentation \eqref{eq_presentation-M}
induces a presentation of sets :
$$
\Pol^{\underline d}_A(\underline M,Q)\isom
\xymatrix{
\Ker\,(\Pol^{\underline d}_A(\underline L,Q) \ar@<.5ex>[r]
\ar@<-.5ex>[r]& \Pol^{\underline d}_A(\underline L',Q))
}\qquad
\text{for every $A$-module $Q$}.
$$
\end{claim}
\begin{pfclaim} Indeed, suppose that
$\lambda:\underline L\squig Q$ is a given homogeneous
multipolynomial law of degree $\underline d$, such that
$\lambda\circ\underline h^1=\lambda\circ\underline h^2$
(notation of remark \ref{rem_apriori}(i)).
Since the induced map
$L'_i\otimes_AB\to(L_i\otimes_AB)\times_{M\otimes_AB}(L_i\otimes_AB)$
is surjective for every $i=1,\dots,k$, it follows that
$\lambda_B$ descends (uniquely) to a map
$\mu_B:(M_1\times\cdots\times M_k)\otimes_AB\to Q\otimes_AB$.
An easy inspection shows that the resulting rule :
$B\mapsto\mu_B$ is a homogeneous multipolynomial law of
degree $\underline d$.
\end{pfclaim}

From claim \ref{cl_presentation-M} and lemma
\ref{lem_multipolynomial}, it follows easily that
the theorem holds for $\underline M$, with
$$
\xymatrix{
\Gamma^{\underline d}_A(\underline M):=\Coker\,(
\Gamma^{\underline d}_A(\underline h^1),
\Gamma^{\underline d}_A(\underline h^2):
\Gamma^{\underline d}_A(\underline L')
\ar@<.5ex>[r] \ar@<-.5ex>[r] & \Gamma^{\underline d}_A(\underline L))
}$$
(details left to the reader). Lastly, let $\underline M$
be an arbitrary sequence. We may then find a filtered
system of sequences $(\underline M{}_\sigma~|~\sigma\in\Sigma)$
consisting of finitely presented $A$-modules and with
$A$-linear transition maps, whose colimit is isomorphic
to $\underline M$. Since such colimits are preserved by
arbitrary base changes, and commute with the forgetful
functor to sets, it follows easily that the system induces
an isomorphism of sets :
$$
\Pol^{\underline d}_A(\underline M,Q)\isom\lim_{\sigma\in\Sigma}
\Pol^{\underline d}_A(\underline M{}_\sigma,Q)
\qquad
\text{for every $A$-module $Q$}
$$
(details left to the reader). In view of the foregoing, it
then follows that the theorem holds for $\underline M$, with
$$
\Gamma^{\underline d}_A(\underline M):=\colim_{\sigma\in\Sigma}
\Gamma^{\underline d}_A(\underline M{}_\sigma).
$$
Again, we invite the reader to spell out the details.
\end{proof}

\sset\subsubsection{}\label{subsec_comultipl-Gamma}
The identity map of $\Gamma^{\underline d}_A\underline M$
corresponds to a homogeneous multipolynomial law
$$
\blambda_{\underline M}^{\underline d}:
\underline M\squig\Gamma^{\underline d}_A\underline M
$$
such that every other law $\underline M\squig Q$ homogeneous
of the same degree, factors uniquely through
$\blambda^{\underline d}_{\underline M}$ and an $A$-linear map
$\Gamma^{\underline d}_A\underline M\to Q$ (in the sense
explained in remark \ref{rem_apriori}(i)). Especially,
if $\underline M'\in\Ob((A\Mod)^{k'})$ is another sequence,
and $\underline d'\in\N^{\oplus k'}$ any sequence of integers,
the tensor product $\blambda^{\underline d}_{\underline M}\otimes 
\blambda^{\underline d'}_{\underline M'}$ (remark
\ref{rem_apriori}(ii)) factors uniquely through
$\blambda^{\underline d,\underline d'}_{\underline M,\underline M'}$
and an $A$-linear map
\set\begin{equation}\label{eq_comulti-Gamma}
\Gamma_A^{\underline d,\underline d'}(\underline M,\underline M')
\to\Gamma_A^{\underline d}(\underline M)\otimes_A
\Gamma_A^{\underline d'}(\underline M').
\end{equation}
Furthermore, notice that, in the situation of remark
\ref{rem_apriori}(iii), the mapping \eqref{eq_base-change-pol}
is clearly $A$-linear; hence, we get an induced $B$-linear map
\set\begin{equation}\label{eq_base-change-Gamma}
B\otimes_A\Gamma^{\underline d}_A(\underline M)\to
\Gamma^{\underline d}_B(B\otimes_A\underline M)
\end{equation}
for every $\underline M$ and $\underline d$ as above, and
every $A$-algebra $B$.

\begin{corollary}\label{cor_gammas-multi}
With the notation of \eqref{subsec_comultipl-Gamma}, we have :
\begin{enumerate}
\item
The functor \eqref{eq_Gamma-multi} commutes with all
filtered colimits.
\item
The map \eqref{eq_comulti-Gamma} is an isomorphism of
$A$-modules, for every $\underline M$, $\underline M'$,
$\underline d$ and $\underline d'$.
\item
For every $A$-algebra $B$, the map \eqref{eq_base-change-Gamma}
is an isomorphism of $B$-modules.
\item
If the $A$-modules $M_1,\dots,M_k$ are all free (resp.
finitely presented, resp. finitely generated, resp. flat),
then the same holds for\/ $\Gamma^{\underline d}_A(\underline M)$.
\end{enumerate}
\end{corollary}
\begin{proof}(i): It suffices to check that the natural map
$$
\Pol^{\underline d}_A(\colim_{\sigma\in\Sigma}\underline M{}_\sigma,Q)
\to\lim_{\sigma\in\Sigma}\Pol^{\underline d}_A(\underline M{}_\sigma,Q)
$$
is an isomorphism, for every filtered system of $A$-modules
$(\underline M{}_\sigma~|~\sigma\in\Sigma)$ and every $A$-module
$Q$; the latter assertion is obvious (details left to the reader).

(iv): The assertion concerning free modules of finite rank,
and finitely generated or finitely presented modules follows
directly from the explicit construction in lemma
\ref{lem_multipolynomial} and theorem \ref{th_construct-Gammas}.
The assertion for flat modules follows from the assertion
for free modules of finite rank and from (i), by means of
Lazard's theorem \cite[Ch.I, Th.1.2]{La}. Moreover, from
remark \ref{rem_funct-on-M}(iii), we see that if
$\phi:\underline M\to\underline N$ is a morphism of sequences
of free $A$-modules of finite rank, such that $\phi_i$ is
a split injection for every $i=1,\dots,n$, then
$\Gamma^{\underline d}_A(\phi)$ is also a split injective map
of free $A$-modules; from this, and from (i), it follows
easily that $\Gamma^{\underline d}_A$ transforms sequences
of free $A$-modules (of possibly infinite rank) into
free $A$-modules : details left to the reader.

(iii): In light of (i), it is easily seen that the
natural transformation \eqref{eq_comulti-Gamma} commutes
with arbitrary filtered colimits of $A$-modules, hence
we are reduced to the case where $\underline M$ is a sequence
of finitely presented $A$-modules. Next, claim
\ref{cl_presentation-M} shows that a presentation
\eqref{eq_presentation-M} yields a commutative diagram
$$
\xymatrix{
\Gamma^{\underline d}_A(\underline L')
\ar@<.5ex>[r] \ar@<-.5ex>[r] \ar[d] &
\Gamma^{\underline d}_A(\underline L) \ar[r] \ar[d] &
\Gamma^{\underline d}_A(\underline M) \ar[r] \ar[d] & 0 \\
\Gamma^{\underline d}_B(B\otimes_A\underline L')
\ar@<.5ex>[r] \ar@<-.5ex>[r] &
\Gamma^{\underline d}_B(B\otimes_A\underline L) \ar[r] &
\Gamma^{\underline d}_B(B\otimes_A\underline M) \ar[r] & 0
}$$
with exact rows. Therefore, we are further reduced to
the case where $\underline M$ is a sequence of free
$A$-modules of finite rank, in which case the assertion
follows by a simple inspection of the proof of lemma
\ref{lem_multipolynomial}.

(ii): Let $t:
\Gamma^{\underline d,\underline d'}_A(\underline M,\underline M')
\to Q$ be an $A$-linear map. To $t$ we attach an $A$-linear
map
$$
t^*:\Gamma^{\underline d}_A(\underline M)\otimes_A
\Gamma^{\underline d'}_A(\underline M')\to Q
$$
as follows. Set $\tau:=t\circ
\blambda^{\underline d,\underline d'}_{\underline M,\underline M'}:
(\underline M,\underline M')\squig Q$ (notation of remark
\ref{rem_apriori}(i)). To every $A$-algebra $B$, the law
$\tau$ attaches the mapping
$$
B\otimes_A\underline M\to
\Pol^{\underline d'}_B(B\otimes_A\underline M',B\otimes_AQ)
\qquad
x\mapsto\tau_{B,x}
$$
where $(\tau_{B,x})_C(y):=\tau_C(1\otimes x,y)$ for every
$x\in B\otimes_A\underline M$, every $B$-algebra $C$
and every $y\in C\otimes_A\underline M'$. In light of
(iii), the latter is the same as a mapping
$$
B\otimes_A\underline M\to
\Hom_B(B\otimes_A\Gamma^{\underline d'}_A(\underline M'),B\otimes_AQ)
\qquad
x\mapsto t_{B,x}
$$
and notice that, for every $x\in B\otimes_A\underline M$, the
$B$-linear map $t_{B,x}$ is characterized by the identity
\set\begin{equation}\label{eq_charact-one}
t_{B,x}(\blambda^{\underline d'}_{\underline M',B}(y))=
(\tau_{B,x})_B(y)
\qquad
\text{for every $y\in B\otimes_A\underline M'$}.
\end{equation}
We deduce a compatible system of mappings
$$
B\otimes_A\Gamma^{\underline d'}_A(\underline M')\to
\Pol_B^{\underline d}(B\otimes_A\underline M,B\otimes_AQ)
\quad :\quad
y\mapsto\tau^*_{B,y}
\qquad
\text{for every $A$-algebra $B$}
$$
where $(\tau^*_{B,y})_C(x):=t_{C,x}(1\otimes y)$ for every
$B$-algebra $C$, every $x\in C\otimes_A\underline M$, and every
$y\in B\otimes_A\Gamma^{\underline d'}_A(\underline M)$. The
$C$-linearity of $t_{C,x}$ implies that these latter mappings
are $B$-linear (for every $A$-algebra $B$), and then they
are the same as a compatible system of $B$-linear maps
\set\begin{equation}\label{eq_same-as-t*}
B\otimes_A\Gamma^{\underline d'}_A(\underline M')\to
\Hom_B(B\otimes_A\Gamma^{\underline d}_A(\underline M),B\otimes_AQ)
\quad :\quad
y\mapsto t^*_{B,y}
\end{equation}
for every $A$-algebra $B$. Notice that, for every
$y\in B\otimes_A\Gamma^{\underline d'}_A(\underline M')$,
the $B$-linear map $t^*_{B,y}$ is characterized by the identity :
\set\begin{equation}\label{eq_charact-two}
t^*_{B,y}(\blambda^{\underline d}_{\underline M,B}(x))=
(\tau^*_{B,y})_B(x)
\qquad
\text{for every $x\in B\otimes_A\underline M$}.
\end{equation}
Obviously, the datum of a compatible system \eqref{eq_same-as-t*}
is the same as that of a map $t^*$ as sought.

\begin{claim}\label{cl_t-and-taus}
$t^*\circ\eqref{eq_comulti-Gamma}=t$.
\end{claim}
\begin{pfclaim} It suffices to show that
$t^*\circ\eqref{eq_comulti-Gamma}\circ
\blambda^{\underline d,\underline d'}_{\underline M,\underline M'}=
\tau$. However, recall that $\eqref{eq_comulti-Gamma}\circ
\blambda^{\underline d,\underline d'}_{\underline M,\underline M'}=
\blambda^{\underline d}_{\underline M}\otimes
\blambda^{\underline d'}_{\underline M'}$, so we are reduced to
checking that
$t^*\circ(\blambda^{\underline d}_{\underline M}\otimes
\blambda^{\underline d'}_{\underline M'})=\tau$.
Hence, let $B$ be any $A$-algebra, and
$m\in B\otimes_A\underline M$, $m'\in B\otimes_A\underline M'$
any two elements. To ease notation, set
$x:=\blambda^{\underline d}_{\underline M}(m)$ and
$y:=\blambda^{\underline d'}_{\underline M'}(m')$. We compute :
$$
\begin{aligned}
(t^*\circ(\blambda^{\underline d}_{\underline M}\otimes
\blambda^{\underline d'}_{\underline M'}))_B(m,m')=\: &
(B\otimes_At^*)(x\otimes y) \\
=\: & t^*_{B,y}(x) \\
=\: & (\tau^*_{B,y})_B(m) & & \text{(by \eqref{eq_charact-two})} \\
=\: & t_{B,m}(y) \\
=\: & (\tau_{B,m})_B(m') & & \text{(by \eqref{eq_charact-one})} \\
=\: & \tau_B(m,m')
\end{aligned}
$$
as claimed.
\end{pfclaim}

Especially, if we let $t$ be the identity map of
$\Gamma_A^{\underline d,\underline d'}(\underline M,\underline M')$,
claim \ref{cl_t-and-taus} shows that the resulting $t^*$ is a left
inverse for \eqref{eq_comulti-Gamma}. Now, if all the modules
$M_1,\dots,M'_{k'}$ are free of finite rank, then lemma
\ref{lem_multipolynomial} shows that both of the $A$-modules
appearing in \eqref{eq_comulti-Gamma} are free and have the
same rank, so the assertion follows, in this case.

Next, suppose that the modules of the two sequences are
finitely presented; in this case, we argue by descending
induction on the number $f$ of free $A$-modules appearing
in the two sequences. If $f=k+k'$, the assertion has just
been shown. Suppose that $f<k+k'$, and the assertion
is already known for all pairs of sequences containing
at least $f+1$ free modules. After permutation of the
sequences, we may assume that $M_i$ is not free, for
some index $i\leq k$. Then, pick a presentation of
$M_i$ as in \eqref{eq_presentation-M}, and let
$\underline N'$ (resp. $\underline N''$) be the
sequences obtained from $\underline M$ by replacing
$M_i$ with $L_i$ (resp. with $L'_i$); by naturality of
\eqref{eq_comulti-Gamma}, the maps $h^1_i$
and $h^2_i$ induce a commutative diagram with exact rows :
$$
\xymatrix{
\Gamma^{\underline d}_A(\underline N'',\underline M')
\ar@<.5ex>[r] \ar@<-.5ex>[r] \ar[d] &
\Gamma^{\underline d}_A(\underline N',\underline M') \ar[r] \ar[d] &
\Gamma^{\underline d}_A(\underline M,\underline M') \ar[r] \ar[d] & 0 \\
\Gamma^{\underline d}_A(\underline N'')\otimes_A
\Gamma^{\underline d}_A(\underline M')
\ar@<.5ex>[r] \ar@<-.5ex>[r] &
\Gamma^{\underline d}_A(\underline N')\otimes_A
\Gamma^{\underline d}_A(\underline M') \ar[r] &
\Gamma^{\underline d}_A(\underline M)\otimes_A
\Gamma^{\underline d}_A(\underline M') \ar[r] & 0
}$$
so it suffices to prove the assertion for the pair
of sequences $(N',M')$ and $(N'',M')$; but the latter
contain $f+1$ free $A$-modules, so we conclude
by induction.

Lastly, if $\underline M$ and $\underline N$ are arbitrary
sequences, we may present them as filtered colimits of
sequences of finitely presented $A$-modules; since
tensor products commute with such colimits, we are done,
by the foregoing case (details left to the reader).
\end{proof}

\begin{remark}\label{rem_compute-Gamma-Pol}
(i)\ \ 
Let $\underline M$ and $\underline d$ be as in
\eqref{subsec_comultipl-Gamma}, suppose that all the
$A$-modules $M_i$ are free of finite rank, and set
$S(\underline M):=\Sym^\bullet_A(\underline M^\vee)$.
Let $Q$ be another $A$-module, and
$\lambda:\underline M\squig Q$ a homogeneous multipolynomial
law of degree $\underline d$. From lemma
\ref{lem_multipolynomial} we may extract the following
calculation for the $A$-linear map
$\lambda^\dagger:\Gamma^{\underline d}_A(\underline M)\to Q$
corresponding to $\lambda$.
Set
$$
P_\lambda:=
\lambda_{S(\underline M)}(\one_{M_1},\dots,\one_{M_k})
\in\Gamma^{\underline d}_A(\underline M)^\vee\otimes_AQ.
$$
Then $\lambda^\dagger$ is the image of $P_\lambda$ under the
natural identification
$\Gamma^{\underline d}_A(\underline M)^\vee\otimes_AQ\isom
\Hom_A(\Gamma^{\underline d}_A(\underline M),Q)$.

(ii)\ \
Let $\underline M$, $\underline M'$, $\underline d$ and
$\underline d'$ be as in \eqref{subsec_comultipl-Gamma},
suppose that all the $A$-modules $M_i$ and $M'_i$ are
free of finite rank, and to ease notation let
$$
\blambda:=\blambda^{\underline d}_{\underline M}
\qquad
\blambda':=\blambda^{\underline d'}_{\underline M'}
\qquad
\lambda:=\blambda\otimes\blambda'.
$$
By (i), we compute
$$
P_\lambda=
\blambda_{S(\underline M,\underline M')}(\one_{M_1},\dots,\one_{M_k})
\otimes
\blambda'_{S(\underline M,\underline M')}(\one_{M'_1},\dots,\one_{M'_{k'}}).
$$
However,
$\blambda_{S(\underline M,\underline M')}(\one_{M_1},\dots,\one_{M_k})$
is the image of $P_\blambda$ under the $A$-linear induced by the
$A$-algebra homomorphism
$S(\underline M)\to S(\underline M,\underline M')=
S(\underline M)\otimes_AS(\underline M')$ given by the rule
$s\mapsto s\otimes 1$, for every $s\in S(\underline M)$.
Likewise we determine
$\blambda'_{S(\underline M,\underline M')}(\one_{M'_1},\dots,\one_{M'_{k'}})$,
and we conclude that $P_\lambda=P_\blambda\otimes P_{\blambda'}$,
whence $\lambda^\dagger=\blambda^\dagger\otimes\blambda^\dagger$,
under the natural identification
$$
\End_A(\Gamma^{\underline d,\underline d'}_A(\underline M,\underline M'))
\isom
\End_A(\Gamma^{\underline d}_A(\underline M))\otimes_A
\End_A(\Gamma^{\underline d'}_A(\underline M')).
$$
But by definition, $(\blambda^{\underline d}_{\underline M})^\dagger$
is the identity automorphism of $\Gamma^{\underline d}_A(\underline M)$,
for every sequence $\underline M$, and every degree $\underline d$.
We see therefore that :
$$
(\blambda^{\underline d}_{\underline M}\otimes
\blambda^{\underline d'}_{\underline M'})^\dagger=
(\blambda^{\underline d,\underline d'}_{\underline M,\underline M'})^\dagger
$$
which translates as saying that the $A$-linear map
\eqref{eq_comulti-Gamma} is none else than the transpose
of the natural isomorphism
$$
\Sym^{\underline d}_A(\underline M^\vee)\otimes_A
\Sym^{\underline d'}_A(\underline M'{}^\vee)\isom
\Sym^{\underline d,\underline d'}_A(\underline M^\vee,\underline M'{}^\vee).
$$
This gives an alternative way to prove corollary
\ref{cor_gammas-multi}(ii), in the case of free $A$-modules.
\end{remark}

\sset\subsubsection{}\label{subsec_late-insight}
Let $\underline M$ and $\underline M'$ be two sequences
of $A$-modules, of the same length $k$, and $\underline e$,
$\underline e'$ two degrees of length $k$. If
$\lambda:(\underline M,\underline M')\squig Q$ is any
homogeneous law of degree $(\underline e,\underline e')$,
then we can also regard $\lambda$ as a homogeneous
multipolynomial law of degree $\underline e+\underline e'$
on the sequence $\underline M\oplus\underline M'$ of
length $k$. In this way, we obtain a natural $A$-linear mapping
\set\begin{equation}\label{eq_need-this-too}
\Pol_A^{\underline e,\underline e'}((\underline M,\underline M'),Q)
\to\Pol^{\underline e+\underline e'}_A(\underline M\oplus\underline M',Q)
\qquad
\text{for every $A$-module $Q$}
\end{equation}
which, by transposition, corresponds to a unique $A$-linear map
\set\begin{equation}\label{eq_saturday}
\Gamma^{\underline e+\underline e'}_A(\underline M\oplus\underline M')
\to\Gamma^{\underline e,\underline e'}_A(\underline M,\underline M').
\end{equation}
Fix $\underline d\in\N^{\oplus k}$; in view of corollary
\ref{cor_gammas-multi}(ii), we see that the sum of the
maps \eqref{eq_saturday} for all pairs $(\underline e,\underline e')$
with $\underline e+\underline e'=\underline d$, is a well
defined $A$-linear map
$$
\Delta^{\underline d}_{\underline M,\underline M'}:=
\bigoplus_{\underline e+\underline e'=\underline d}
\Delta^{\underline e,\underline e'}_{\underline M,\underline M'}:
\Gamma^{\underline d}_A(\underline M\oplus\underline M')\to
\bigoplus_{\underline e+\underline e'=\underline d}
\Gamma^{\underline e}_A(\underline M)\otimes_A
\Gamma^{\underline e'}_A(\underline M').
$$

\begin{proposition}\label{prop_Deltas-are-iso}
The map $\Delta^{\underline d}_{\underline M,\underline M'}$ is an
isomorphism, for every $\underline M$, $\underline M'$,
and every degree $\underline d$.
\end{proposition}
\begin{proof} Arguing as in the proof of corollary
\ref{cor_gammas-multi}(ii), we are easily reduced to the
case where $\underline M$ and $\underline M'$ are sequences
of free $A$-modules of finite rank. In this case,
\eqref{eq_saturday} can be computed as in remark
\ref{rem_compute-Gamma-Pol}(i); namely, it corresponds
to the element
$$
P\in
\Sym^{\underline e+\underline e'}_A((\underline M\oplus\underline M')^\vee)
\otimes_A
\Sym_A^{\underline e,\underline e'}(\underline M^\vee,\underline M'{}^\vee)^\vee
$$
determined as follows. Take
$Q:=\Gamma^{\underline e,\underline e'}_A(\underline M,\underline M')$,
and let
$\lambda:\underline M\oplus\underline M'\squig
\Gamma^{\underline e,\underline e'}_A(\underline M,\underline M')$
be the image of
$\blambda^{\underline e,\underline e'}_{\underline M,\underline M'}$
under \eqref{eq_need-this-too}; then
$$
P=\lambda_{S(\underline M\oplus\underline M')}
(\one_{M_1\oplus M'_1},\dots,\one_{M_k\oplus M'_k}).
$$
Thus, we come down to evaluating the image of $\one_{M_i\oplus M'_i}$
under the natural identification
$$
(M_i\oplus M'_i)\otimes_AS(M_i\oplus M'_i)
\isom
(M_i\otimes_AS(M_i\oplus M'_i))\times(M'_i\otimes_AS(M_i\oplus M'_i)).
$$
for every $i=1,\dots,k$. We easily find that
$\one_{M_i\oplus M'_i}\mapsto(p_i,p'_i)$ under
this identification, where $p_i:M_i\oplus M'_i\to M_i$
and $p'_i:M_i\oplus M'_i\to M'_i$ are the natural
projections; so finally
\set\begin{equation}\label{eq_uffa}
P=(\blambda^{\underline e,\underline e'}_{\underline M,
\underline M'})_{S(\underline M\oplus\underline M')}
(p_1,\dots,p_k,p'_1,\dots,p'_k).
\end{equation}
Pick a basis $(e_{i1},\dots,e_{ir_i})$ (resp.
$(e'_{i1},\dots,e'_{ir'_i})$) for each $M_i$ (resp. $M'_i$),
let $(e^*_{i1},\dots,e^*_{ir_i})$ (resp.
$(e^{\prime*}_{i1},\dots,e^{\prime*}_{ir'_i})$)) be the dual basis,
and for every $=1,\dots,k$ and every $j=1,\dots,r_i$
(resp. $j=1,\dots,r'_i$) denote by $\eps_{ij}$ (resp.
$\eps'_{ij}$) the element of $(M_i\oplus M'_i)^\vee$
that agrees with $e^*_{ij}$ on $M_i$ and that vanishes
on $M'_i$ (resp. that agrees with $e^{\prime*}_{ij}$ on
$M'_i$ and vanishes on $M_i$). Notice that
$$
p_i=\sum^{r_i}_{j=1}\eps^*_{ij}\otimes e_{ij}
\qquad
p'_i=\sum^{r'_i}_{j=1}\eps^{\prime*}_{ij}\otimes e'_{ij}
\qquad
\text{for every $i=1,\dots,k$}.
$$
With this notation, it follows that the right hand-side
of \eqref{eq_uffa} is the image of
$$
(\blambda^{\underline e,\underline e'}_{\underline M,
\underline M'})_{S(\underline M,\underline M')}(\one_{M_1},\dots,\one_{M'_k})
=\one_{\Sym^{\underline e,\underline e'}_A(\underline M,\underline M')}
$$
under the natural map
\set\begin{equation}\label{eq_see_clear_now}
S(\underline M,\underline M')\to S(\underline M\oplus\underline M')
\quad :\quad
e^*_{ij}\mapsto\eps^*_{ij}
\quad
e^{\prime*}_{ij}\mapsto\eps^{\prime*}_{ij}.
\end{equation}
In other words, the sought $P$ is none else than the
restriction of \eqref{eq_see_clear_now} to the direct
factor $\Sym_A^{\underline e,\underline e'}
(\underline M^\vee,\underline M'{}^\vee)$, and this
restriction is an isomorphism onto a direct factor of
the $A$-module
$\Sym^{\underline e+\underline e'}_A(\underline M\oplus\underline M')$.
By transposition, we see that \eqref{eq_saturday} is identified
with the projection onto a natural direct factor of
$\Gamma^{\underline e+\underline e'}_A(\underline M\oplus\underline M')$,
and a simple inspection shows that the sum of all these projections
is the identity map of the latter $A$-module, whence the proposition.
\end{proof}

\begin{remark}\label{rem_Groth-style}
(i)\ \ 
Let $\underline M,\underline M'$, $Q$ and
$\lambda$ be as in \eqref{subsec_late-insight}. The image
of $\lambda$ under \eqref{eq_need-this-too} is characterized
as the unique law $\phi:\underline M\oplus\underline M'\squig Q$
of degree $\underline e+\underline e'$ such that
$\phi(m,m')=\lambda_B(m,m')$ for every $A$-algebra $B$ and
every $(m,m')\in B\otimes_A(\underline M\oplus\underline M')$.
In turns, $\phi$ corresponds to the unique $A$-linear map
$f:\Gamma^{\underline e+\underline e'}_A(\underline M\oplus\underline M')
\to Q$ such that 
$(B\otimes_Af)(
\blambda^{\underline e+\underline e'}_{\underline M\oplus\underline M'})_B
(m,m')=\lambda_B(m,m')$ for every $B$ and $(m,m')$ as above.
Especially, if we take $\lambda:=
\lambda^{\underline e,\underline e'}_{\underline M,\underline M'}$,
we see that \eqref{eq_saturday} is characterized as the unique
$A$-linear map such that
$$
(\blambda^{\underline e+\underline e'}_{\underline M\oplus\underline M'})_B
(m,m')\mapsto
(\blambda^{\underline e,\underline e'}_{\underline M,\underline M'})_B(m,m')
$$
for every $B$ and $(m,m')$ as above. Lastly, it follows
that $\Delta^{\underline e,\underline e'}_{\underline M,\underline M'}$
is characterized as the unique $A$-linear map such that
$$
(\blambda^{\underline e+\underline e'}_{\underline M\oplus\underline M'})_B
(m,m')\mapsto(\blambda^{\underline e}_{\underline M})_B(m)\otimes
(\lambda^{\underline e'}_{\underline M'})_B(m')
$$
for every $A$-algebra $B$ and every
$(m,m')\in B\otimes_A(\underline M\oplus\underline M')$.

(ii)\ \
In the same vein, let $f:\underline M\to\underline N$ be
any morphism in $(A\Mod)^a$, and $\underline d\in\N^{\oplus k}$
any degree. Then we see that the induced map
$$
\Gamma^{\underline d}_A(f):\Gamma^{\underline d}_A(\underline M)
\to\Gamma^{\underline d}_A(\underline N)
$$
is characterized as the unique $A$-linear map such that
$$
(\blambda^{\underline d}_{\underline M})_B(m)\mapsto
(\blambda^{\underline d}_{\underline N})_B(f(m))
$$
for every $A$-algebra $B$, and every $m\in B\otimes_A\underline M$
(details left to the reader).
\end{remark}

\sset\subsubsection{}\label{subsec_for-equation}
Let $\underline M$, $\underline M'$, $\underline M''$ be
three sequences of $A$-modules of length $k$, and $\underline d$,
$\underline d'$, $\underline d''$ three degrees, also of lengths
$k$. We obtain a diagram of $A$-linear maps :
$$
\xymatrix{ \Gamma_A^{\underline e+\underline e'+\underline e''}
(\underline M\oplus\underline M'\oplus\underline M'')
\ar[rrrr]^-{\Delta^{\underline e,\underline e'+
\underline e''}_{\underline M,\underline M'\oplus\underline M''}}
\ar[d]_{\Delta^{\underline e+\underline e',
\underline e''}_{\underline M\oplus \underline M',\underline M''}} & & & &
\Gamma_A^{\underline e}(\underline M)\otimes_A
\Gamma^{\underline e'+\underline e''}_A
(\underline M'\oplus\underline M'')
\ar[d]^{\Gamma^{\underline e}_A(\underline M)\otimes_A
\Delta^{\underline e',\underline e''}_{\underline M',\underline M''}} \\
\Gamma^{\underline e+\underline e'}_A(\underline M\oplus\underline M')
\otimes_A\Gamma^{\underline e''}_A(\underline M'')
\ar[rrrr]^-{\Delta^{\underline e,\underline e'}_{\underline M,
\underline M'}\otimes_A\Gamma^{\underline e''}_A(\underline M'')} & & & &
\Gamma^{\underline e}_A(\underline M)\otimes_A
\Gamma^{\underline e'}_A(\underline M')\otimes_A
\Gamma^{\underline e''}_A(\underline M'').
}$$

\begin{proposition}
The diagram of \eqref{subsec_for-equation} commutes.
\end{proposition}
\begin{proof} By remark \ref{rem_Groth-style}(i), the map
$\Delta^{\underline e,\underline e'+
\underline e''}_{\underline M,\underline M'\oplus\underline M''}$
is characterized as the unique one such that
$$
(\blambda^{\underline e+\underline e'+\underline e''}_{\underline M
\oplus\underline M'})_B(m,m',m'')\mapsto
(\blambda^{\underline e}_{\underline M})_B(m)\otimes
(\blambda^{\underline e'+\underline e''}_{\underline M\oplus\underline M'})_B(m',m'')
$$
for every $A$-algebra $B$ and every $(m,m',m'')\in B\otimes_A
(\underline M\oplus\underline M'\oplus\underline M'')$, and
$\Delta^{\underline e',\underline e''}_{\underline M',\underline M''}$
is the unique $A$-linear map such that
$$
(\blambda^{\underline e'+\underline e''}_{\underline M'\oplus\underline M''})_B
(m',m'')\mapsto(\blambda^{\underline e'}_{\underline M'})_B(m')\otimes
(\lambda^{\underline e''}_{\underline M''})_B(m'')
$$
for every $A$-algebra $B$ and every
$(m',m'')\in B\otimes_A(\underline M'\oplus\underline M'')$.
Therefore, the composition of the top horizontal and right
vertical arrows in the diagram is the unique $A$-linear map
such that
$$
(\blambda^{\underline e+\underline e'+\underline e''}_{\underline M
\oplus\underline M'})_B(m,m',m'')\mapsto
(\blambda^{\underline e}_{\underline M})_B(m)\otimes
(\blambda^{\underline e'}_{\underline M'})_B(m')\otimes
(\lambda^{\underline e''}_{\underline M''})_B(m'')
$$
for every $B$ and $(m,m',m'')$ as above.
The same calculation can be carried out for the composition
of the other two arrows, and clearly one obtains the same
result.
\end{proof}

\sset\subsubsection{}\label{subsec_divided-power-defined}
Let $\underline M$ be any sequence of $A$-modules of length
$k$, and denote by $s:\underline M\oplus\underline M\to\underline M$
the addition map : $(m,m')\mapsto m+m'$ for every
$m,m'\in\underline M$. The map $s$ induces $A$-linear maps
$$
\mu_{\underline M}^{\underline e,\underline e'}:
\Gamma^{\underline e}_A(\underline M)\otimes_A
\Gamma^{\underline e'}(\underline M)
\xrightarrow{\ \delta^{\underline e,\underline e'}\ }
\Gamma^{\underline e+\underline e'}_A(\underline M\oplus\underline M)
\xrightarrow{\ \Gamma^{\underline e+\underline e'}_A(s)\ }
\Gamma^{\underline d}_A(\underline M)
\qquad
\text{for every $\underline e,\underline e'\in\N^{\oplus k}$}
$$
where $\delta^{\underline e,\underline e'}$ is the natural
inclusion map obtained by restricting the inverse of
$\Delta^{\underline e+\underline e'}_{\underline M,\underline M}$
to the direct factor $\Gamma^{\underline e}_A(\underline M)
\otimes_A\Gamma^{\underline e'}(\underline M)$.
In light of remark \ref{rem_Groth-style}(i,ii), the map
$\mu_{\underline M}^{\underline e,\underline e'}$ is characterized
as the unique $A$-linear map such that
\set\begin{equation}\label{eq_almost-divided-powers}
(\blambda^{\underline e}_{\underline M})_B(m)\otimes
(\blambda^{\underline e'}_{\underline M})_B(m')\mapsto
(\blambda^{\underline e+\underline e'}_{\underline M})_B(m+m')
\end{equation}
for every $A$-algebra $B$ and every
$(m,m')\in B\otimes_A(\underline M\oplus\underline M)$.
Using this characterization, it is easily seen that the
diagram
$$
\xymatrix{ \Gamma^{\underline e}_A(\underline M)\otimes_A
\Gamma^{\underline e'}_A(\underline M)\otimes_A
\Gamma^{\underline e''}_A(\underline M)
\ar[rrr]^-{\mu^{\underline e,\underline e'}_{\underline M}
\otimes_A\Gamma^{\underline e''}_A(\underline M)}
\ar[d]_{\Gamma^{\underline e}_A(\underline M)\otimes_A
\mu^{\underline e',\underline e''}_{\underline M}} & & &
\Gamma^{\underline e+\underline e'}_A(\underline M)\otimes_A
\Gamma^{\underline e''}_A(\underline M)
\ar[d]^{\mu^{\underline e+\underline e',\underline e''}_{\underline M}} \\
\Gamma^{\underline e}_A(\underline M)\otimes_A
\Gamma^{\underline e'+\underline e''}_A(\underline M)
\ar[rrr]^-{\mu^{\underline e,\underline e'+\underline e''}_{\underline M}}
& & & \Gamma^{\underline e+\underline e'+\underline e''}_A(\underline M)
}$$
commutes, for every
$\underline e,\underline e',\underline e''\in\N^{\oplus k}$.
The foregoing shows that the maps $\mu^{\bullet\bullet}_{\underline M}$
define on
$$
\Gamma^\bullet_A(\underline M):=\bigoplus_{\underline d\in\N^{\oplus k}}
\Gamma^{\underline d}_A(\underline M).
$$
a natural structure of associative $\N^{\oplus k}$-graded $A$-algebra.
Moreover, \eqref{eq_almost-divided-powers} also easily implies that
this algebra is commutative. We call $\Gamma^\bullet_A(\underline M)$
the {\em divided power envelope} of $\underline M$. It is
customary to use the notation
$$
m^{[\underline d]}:=(\blambda_{\underline M}^{\underline d})_B(m)
\in B\otimes_A\underline M
$$
for every $A$-algebra $B$, every $m\in B\otimes_A\underline M$,
and every $\underline d\in\N^{\oplus k}$. Then, proposition
\ref{prop_Deltas-are-iso} and a simple inspection of the
definitions yield the identity
\set\begin{equation}\label{eq_newton}
(a+b)^{[\underline d]}=\sum_{\underline e+\underline e'=\underline d}
a^{[\underline e]}\cdot b^{[\underline e']}
\end{equation}
for every $A$-algebra $B$ and every
$a,b\in B\otimes_A\underline M$. Clearly, the above construction
yields a well defined functor
\set\begin{equation}\label{eq_div-pow-alg}
(A\Mod)^k\to A\Alg
\quad :\quad
M\mapsto\Gamma_A^\bullet(\underline M)
\end{equation}
and notice that, for every morphism $f:\underline M\to\underline N$
of sequences, we have
$$
(\Gamma^\bullet_Af)(a^{[\underline d]})=(fa)^{[\underline d]}=
(\Gamma^\bullet_Af(a))^{[\underline d]}
$$
for every degree $\underline d$, every $A$-algebra $B$, and
every $a\in B\otimes_A\underline M$.

\begin{theorem}\label{th_divided-powers}
With the notation of \eqref{subsec_divided-power-defined},
we have :
\begin{enumerate}
\item
For every $A$-algebra $B$, there exists a natural
isomorphism of\/ $\N^{\oplus k}$-graded algebras :
$$
B\otimes_A\Gamma^\bullet_A(\underline M)\isom
\Gamma_B^\bullet(B\otimes_AM).
$$
\item
The functor \eqref{eq_div-pow-alg} commutes with filtered colimits.
\item
For any two sequences of $A$-modules $\underline M$ and
$\underline M'$ of the same length, there exists a natural
isomorphism of graded $A$-algebras :
$$
\Gamma^\bullet_A(\underline M\oplus\underline M')\isom
\Gamma_A^\bullet(\underline M)\otimes_A\Gamma^\bullet_A(\underline M')
$$
(with the $\N^{\oplus k}$-grading of the tensor product
defined in the obvious way).
\end{enumerate}
\end{theorem}
\begin{proof}(i) and (ii) follow immediately from corollary
\ref{cor_gammas-multi}(i,iii), and (iii) follows directly
from proposition \ref{prop_Deltas-are-iso}.
\end{proof}

\sset\subsubsection{}\label{subsec_d-factorial}
Consider the special case of a sequence $\underline M$
consisting of a single $A$-module $M$. From \eqref{eq_newton},
by evaluating the expression $(2m)^{[d]}=2^d\cdot m^{[d]}$,
a simple induction on $d$ shows that
\set\begin{equation}\label{eq_d-factorial}
d!\cdot m^{[d]}=m^d
\qquad
\text{in $\Gamma_A^\bullet(M)$, for every $m\in M$}.
\end{equation}

\begin{corollary}\label{cor_divided-powers}
In the situation of \eqref{subsec_d-factorial}, we have
$$
m^{[i]}\cdot m^{[j]}=\binom{i+j}{i}\cdot m^{[i+j]}
\qquad
\text{in $\Gamma^\bullet_A(M)$, for every $m\in M$ and every $i,j\in\N$}.
$$
\end{corollary}
\begin{proof} By considering the $A$-linear map $f:A\to M$
such that $f(1)=m$, we reduce to the case where $M=A$. By
considering the unique ring homomorphism $\Z\to A$ we reduce
further -- in light of theorem \ref{th_divided-powers}(i) --
to the case where $A=\Z$ and $M=\Z$. In this case, $\Gamma_\Z(\Z)$
is a torsion-free $\Z$-module, by lemma \ref{lem_multipolynomial},
so the identity follows easily from \eqref{eq_d-factorial}.
\end{proof}

\subsection{Regular rings}\label{sec_rings-regular}
In a later section, we shall need a criterion for the
regularity of a certain type of extensions of regular
local rings. Such a criterion is stated incorrectly in
\cite[Ch.0, Th.22.5.4]{EGAIV}; our first task is to
supply a corrected version of {\em loc.cit.}

\sset\subsubsection{}\label{subsec_please}
Consider a local ring $A$, with maximal ideal $\fm_A$, and
residue field $k_A:=A/\fm_A$ of characteristic $p>0$. From
\cite[Ch.0, Th.20.5.12(i)]{EGAIV} we get a complex of
$k_A$-vector spaces
\set\begin{equation}\label{eq_finally-exact}
0\to\fm_A/(\fm_A^2+pA)\xrightarrow{\ d_A\ }\Omega^1_{A/\Z}\otimes_Ak_A
\to\Omega^1_{k_A/\Z}\to 0.
\end{equation}

\begin{lemma}\label{lem_finally-exact}
The complex \eqref{eq_finally-exact} is exact.
\end{lemma}
\begin{proof} Set $(A_0,\fm_{A_0}):=(A/pA,\fm_A/pA)$ and
$\F_p:=\Z/p\Z$; the induced sequence of ring homomorphisms
$\F_p\to A_0\to k_A$ yields a distinguished triangle
(\cite[Ch.II, \S2.1.2.1]{Il}) :
$$
\L_{A_0/\F_p}\derotimes_{A_0}k_A\to\L_{k_A/\F_p}\to\L_{k_A/A_0}
\to\L_{A_0/\F_p}\derotimes_{A_0}k_A[1].
$$
However, \cite[Th.6.5.12(ii)]{Ga-Ra} says that $H_i\L_{k_A/\F_p}=0$
for every $i>0$ (one applies {\em loc.cit.} to the valued field
$(k_A,|\cdot|)$, where $|\cdot|$ is the trivial valuation;
more directly, one may observe that $k_A$ is the colimit of a
filtered system of smooth $\F_p$-algebras, since $\F_p$ is
perfect, and then recall that the functor $\L$ commutes with
filtered colimits). We deduce a short exact sequence
$$
0\to H_1\L_{k_A/A_0}\to\Omega^1_{A_0/\F_p}\otimes_{A_0}k_A\to
\Omega^1_{k_A/\F_p}\to 0
$$
which, under the natural identifications
$$
\Omega^1_{k_A/\Z}\isom\Omega^1_{k_A/\F_p}
\qquad
\Omega^1_{A/\Z}\otimes_Ak_A\isom\Omega^1_{A_0/\F_p}\otimes_{A_0}k_A
\qquad
H_1\L_{k_A/A_0}\isom\fm_{A_0}/\fm_{A_0}^2
$$
(\cite[Ch.III, Cor.1.2.8.1]{Il}) becomes \eqref{eq_finally-exact},
up to replacing $d_A$ by $-d_A$ (\cite[Ch.III, Prop.1.2.9]{Il};
details left to the reader).
\end{proof}

\sset\subsubsection{}
Keep the situation of \eqref{subsec_please}, and define
$$
A_2:=A\times_{k_A}W_2(k_A)
$$
where $W_2(k_A)$ is the ring of $2$-truncated Witt vectors of $k_A$,
as in \eqref{subsec_here-we-trunk}, which is augmented over $k_A$,
via the ghost component map $\bar\omega_0:W_2(k_A)\to k_A$.
Notice that
$$
V_1(k_A):=\Ker\,\bar\omega_0
$$
is an ideal of $W_2(k_A)$ such that $V_1(k_A)^2=0$; especially,
it is naturally a $k_A$-vector space, with addition and
scalar multiplication given respectively by the rules :
$$
(0,a)+(0,b):=(0,a+b)
\qquad
x\cdot(0,a):=(x,0)\cdot(0,a)=(0,x^p\cdot a)
$$
for every $a,b,x\in k_A$. In other words, the map
\set\begin{equation}\label{eq_not-perfect}
V_1(k_A)\to k^{1/p}_A
\qquad
(0,a)\mapsto a^{1/p}
\end{equation}
is an isomorphism of $k_A$-vector spaces, and we also see that
$V_1(k_A)$ is a one-dimensional $k_A^{1/p}$-vector space, with
scalar multiplication given by the rule : $a\cdot(0,b):=(0,a^pb)$
for every $a\in k_A^{1/p}$ and $b\in k_A$. By construction, we
have a natural exact sequence of $A_2$-modules :
\set\begin{equation}\label{eq_V_1-W_2}
0\to V_1(k_A)\to A_2\xrightarrow{\ \pi\ } A\to 0.
\end{equation}
Especially, $A_2$ is a local ring, and $\pi$ is a local
ring homomorphism inducing an isomorphism on residue fields.

\sset\subsubsection{}\label{subsec_force-nextpage}
It is easily seen that the rule $(A,\fm_A)\mapsto A_2$ defines
a functor on the category $\mathbf{Local}$ of local rings and
local ring homomorphisms, to the category of (commutative,
unital) rings. Hence, let us set
$$
(A,\fm_A)\mapsto\bar\Omega_A:=\Omega^1_{A_2/\Z}\otimes_{A_2}A.
$$
It follows that the rule $(A,\fm_A)\mapsto(A,\bar\Omega_A)$
yields a functor
$$
\mathbf{Local}\to\AlgMod
$$
to the category whose objects are the pairs $(B,M)$ where
$B$ is a ring, and $M$ is a $B$-module (the morphisms
$(B,M)\to(B',M')$ are the pairs $(\phi,\psi)$ where
$\phi:B\to B'$ is a ring homomorphism, and
$\psi:B'\otimes_BM\to M$ is a $B'$-linear map : cp.
\cite[Def.2.5.22(ii)]{Ga-Ra}).

In view of \eqref{eq_V_1-W_2} and \cite[Ch.0, Th.20.5.12(i)]{EGAIV}
we get an exact sequence of $A$-modules
$$
V_1(k_A)\to\bar\Omega_A\to\Omega^1_{A/\Z}\to 0
$$
whence, after tensoring with $k_A$, a sequence of
$k_A$-vector spaces
\set\begin{equation}\label{eq_forcepage}
0\to V_1(k_A)\xrightarrow{\ j_A\ }\bar\Omega_A\otimes_Ak_A
\xrightarrow{\ \rho\ }\Omega^1_{A/\Z}\otimes_Ak_A\to 0
\end{equation}
which is right exact by construction.

\begin{proposition}\label{prop_forcenext}
With the notation of \eqref{subsec_force-nextpage}, we have :
\begin{enumerate}
\item
If\/ $p\notin\fm_A^2$, then the sequence \eqref{eq_forcepage} is
exact.
\item
If\/ $p\in\fm_A^2$, then $\Ker\,j_A=\{(0,a^p)~|~a\in k_A\}$.
Especially, the isomorphism \eqref{eq_not-perfect} identifies
$\Ker\,j_A$ with the subfield $k_A$ of $k_A^{1/p}$.
\end{enumerate}
\end{proposition}
\begin{proof} In view of lemma \ref{lem_finally-exact}, the
kernel of $\rho$ is naturally identified with the kernel of the
natural map
$$
\fm_{A_2}/(\fm^2_{A_2}+pA_2)\to\fm_A/(\fm_A^2+pA).
$$
On the other hand, it is easily seen that
$\fm_{A_2}=\fm_A\oplus V_1(k_A)$, and under this identification,
the multiplication law of $A_2$ restricts on $\fm_{A_2}$ to the
mapping given by the rule : $(a,b)\cdot(a',b'):=(aa',0)$ for
every $a,a'\in\fm_A$ and $b,b'\in V_1(k_A)$. There follows a
natural isomorphism
\set\begin{equation}\label{eq_under-deco}
\fm_{A_2}/\fm^2_{A_2}\isom(\fm_A/\fm_A^2)\oplus V_1(k_A)
\end{equation}
which identifies the map
$\fm_{A_2}/\fm_{A_2}^2\to\fm_A/\fm^2_A$ deduced from $\pi$
(notation of \eqref{eq_V_1-W_2}), with the projection on the
first factor. Moreover, by inspecting the definitions, we
easily get a commutative diagram
$$
\xymatrix{
V_1(k_A) \ar[r]^-{j_A} \ar[d] & \bar\Omega_A\otimes_Ak_A \\
\fm_{A_2}/\fm^2_{A_2} \ar[r] &
\fm_{A_2}/(\fm^2_{A_2}+pA_2) \ar[u]_{d_{A_2}}
}$$
whose bottom horizotal arrow is the quotient map, and whose
left vertical arrow is the inclusion map of the direct summand
$V_1(k_A)$ resulting from \eqref{eq_under-deco}.

The map $A_2\to\fm_{A_2}/\fm^2_{A_2}$ given by the rule :
$a\mapsto p\bar a:=pa\pmod{\fm^2_{A_2}}$ factors (uniquely)
through a $k_A$-linear map
$$
t_{A_2}:k_A\to\fm_{A_2}/\fm^2_{A_2}
\qquad
a\pmod{\fm_A}\mapsto pa\pmod{\fm_{A_2}^2}
$$
and likewise we may define a map $t_A:k_A\to\fm_A/\fm_A^2$.
The snake lemma then gives an induced map
$\partial:\Ker\,t_A\to V_1(k_A)$, and in view of the foregoing,
the proposition follows from :

\begin{claim} (i)\ \ If $p\notin\fm_A^2$, then $t_A$ is injective.
\begin{enumerate}
\addenu
\item
If $p\in\fm_A^2$, then $\Img\,\partial=\{(0,a^p)~|~a\in k_A\}$.
\end{enumerate}
\end{claim}
\begin{pfclaim}[] (i) is obvious.

(ii): By virtue of \eqref{eq_Verschiebung}, we have $p=(p,(0,1))$
in $A_2$, and if $(a,y)\in A_2$ is any element, then
$p\cdot(a,y)=(pa,(0,a^p))$, so in case $p\in\fm_A^2$, the map
$t_{A_2}$  is given by the rule :
$$
a\mapsto(pa,(0,\bar a{}^p))=(0,(0,\bar a{}^p))
\in(\fm_A/\fm^2_A)\oplus V_1(k_A)
\qquad
\text{for every $a\in A$}.
$$
By the same token, in this case $t_A$ is the zero map. The claim
follows straightforwardly.
\end{pfclaim}
\end{proof}

\sset\subsubsection{}\label{subsec_bOmega}
Keep the notation of \eqref{subsec_force-nextpage}, and
assume that $p\notin\fm_A^2$, so \eqref{eq_forcepage} is a
$k_A$-extension of $\Omega^1_{A/\Z}\otimes_Ak_A$ by $V_1(k_A)$,
by virtue of proposition \ref{prop_forcenext}; hence,
$\eqref{eq_forcepage}\otimes_{k_A}k_A^{1/p}$ is a
$k_A^{1/p}$-extension of the corresponding $k_A^{1/p}$-vector
spaces. Recall now that $V_1(k_A)$ is naturally a
$k_A^{1/p}$-vector space of dimension one, and let
$\psi:V_1(k_A)\otimes_{k_A}k_A^{1/p}\to V_1(k_A)$ be the scalar
multiplication. By push out along $\psi$, we obtain therefore
an extension $\psi*\eqref{eq_forcepage}\otimes_{k_A}k_A^{1/p}$
fitting into a commutative ladder with exact rows :
$$
\xymatrix{
0 \ar[r] & V_1(k_A)\otimes_{k_A}k^{1/p}_A \ar[r] \ar[d]_{\psi_A} &
\bar\Omega_A\otimes_Ak_A^{1/p} \ar[r] \ar[d]_{\bpsi_A} &
\Omega^1_{A/\Z}\otimes_Ak_A^{1/p} \ar[r] \ar[d] & 0 \\
0 \ar[r] & V_1(k_A) \ar[r]^-{\bj_A} & \bOmega_A \ar[r] &
\Omega^1_{A/\Z}\otimes_Ak_A^{1/p} \ar[r] & 0
}$$
(see \cite[\S2.5.5]{Ga-Ra}). We consider now the mapping :
$$
\bd_A:A\to\bOmega_A
\qquad
a\mapsto\bpsi_A(d(a,\tau_{k_A}(\bar a))\otimes 1)
\qquad
\text{for every $a\in A$}
$$
where :
\begin{itemize}
\item
 $\bar a\in k_A$ is the image of $a$ in $k_A$
\item
$\tau_{k_A}$ is the Teichm\"uller mapping (see \eqref{subsec_Teich}),
so $\tau_{k_A}(\bar a)=(\bar a,0)\in W_2(k_A)$
\item
$d:A_2\to\Omega_{A_2/\Z}$ is the universal derivation of $A_2$.
\end{itemize}
Since $\tau_{k_A}$ is multiplicative, the map $\bd_A$ satisfies
Leibniz's rule, {\em i.e.} we have
$$
\bd_A(ab)=\bar a\cdot\bd_A(b)+\bar b\cdot\bd_A(a)
\qquad
\text{for every $a,b\in A$}.
$$
However, $\bd_A$ is not quite a derivation, since additivity
fails for $\tau_{k_A}$, hence also for $\bd_A$. Instead, from
remark \ref{rem_homogeneous-laws}(ii), and recalling that
$(p-1)!=-1$ in $\F_p$, we get the identity :
$$
\begin{aligned}
\tau_{k_A}(\bar a+\bar b)=\: &
\tau_{k_A}(\bar a)+\tau_{k_A}(\bar b)-\sum_{i=1}^{p-1}
\Bigl(0,\frac{\bar a{}^i}{i!}\cdot\frac{\bar b{}^{p-i}}{(p-i)!}\Bigr) \\
=\: &
\tau_{k_A}(\bar a)+\tau_{k_A}(\bar b)-\sum_{i=1}^{p-1}
\frac{\bar a{}^{i/p}}{i!}\cdot\frac{\bar b{}^{1-i/p}}{(p-i)!}\cdot p
\end{aligned}
$$
for every $\bar a,\bar b\in k_A$. On the other hand, notice that
$$
\begin{aligned}
\bpsi_A(d(0,x\cdot p)\otimes 1)=\: &
\bj_A\circ\psi_A(x\cdot p\otimes 1) \\
=\: & x\cdot\bj_A\circ\psi_A(p\otimes 1) \\
=\: & x\cdot\bpsi_A(d(0,p)\otimes 1) \\
=\: & x\cdot\bpsi_A(d(p-(p,0))\otimes 1) \\
=\: & -x\cdot\bd_A(p)
\end{aligned}
$$
for every $x\in k_A^{1/p}$, whence :
$$
\bd_A(a+b)=\bd_A(a)+\bd_A(b)+\sum_{i=1}^{p-1}
\frac{\bar a{}^{i/p}}{i!}\cdot\frac{\bar b{}^{1-i/p}}{(p-i)!}\cdot\bd_A(p)
\qquad
\text{for every $a,b\in A$}.
$$

\sset\subsubsection{}\label{subsec_bar-bd}
Especially, notice we do have $\bd_A(a+b)=\bd_A(a)+\bd_A(b)$ in case
either $a$ or $b$ lies in $\fm_A$. Hence, $\bd_A$ restricts to
an additive map $\fm_A\to\bOmega_A$, and Leibniz's rule implies
that the latter descends to a $k_A$-linear map
$$
\bar\bd_A:\fm_A/\fm_A^2\to\bOmega_A.
$$

\begin{proposition}\label{prop_lift-lemma}
In the situation of \eqref{subsec_please}, suppose that
$p\notin\fm_A^2$. Then there exists a natural exact sequence
of $k^{1/p}_A$-vector spaces :
$$
0\to\fm_A/\fm_A^2\otimes_{k_A}k_A^{1/p}
\xrightarrow{\ \bar\bd_A\otimes k_A^{1/p}\ }\bOmega_A\to
\Omega^1_{k_A/\Z}\otimes_{k_A}k_A^{1/p}\to 0.
$$
\end{proposition}
\begin{proof} By inspecting the constructions in
\eqref{subsec_bOmega}, we obtain a commutative ladder
of $k_A$-vector spaces, with exact rows :
$$
\xymatrix{ 0 \ar[r] & pA/(\fm_A^2\cap pA) \ar[r] \ar[d]_i &
\fm_A/\fm_A^2 \ar[r] \ar[d]^{\bar\bd_A} &
\fm_A/(\fm^2_A+pA) \ar[r] \ar[d]^{d_A} & 0 \\
0 \ar[r] & V_1(k_A) \ar[r]^-{\bj_A} & \bOmega_A \ar[r] &
\Omega^1_{A/\Z}\otimes_Ak_A^{1/p} \ar[r] & 0
}$$
such that $d_A\otimes_{k_A}k_B$ is injective with cokernel
naturally isomorphic to $\Omega^1_{k_A/\Z}\otimes_{k_A}k_A^{1/p}$
(lemma \ref{lem_finally-exact}). By the snake lemma, it
then suffices to show that $i\otimes_{k_A}k_A^{1/p}$ is
an isomorphism; but since both the source and target of
the latter map are one-dimensional $k_A^{1/p}$-vector
spaces, we come down to checking that $i$ is not the zero
map. But a simple inspection yields $i(p)=(0,1)$, as required.
\end{proof}

\sset\subsubsection{}\label{subsec_bOmega-funct}
We wish next to show how the rule $A\mapsto\bOmega_A$
extends to a functor
$$
\mathbf{Local}\to\AlgMod
\qquad
(A,\fm_A)\mapsto(k^{1/p}_A,\bOmega_A).
$$
Namely, if $p$ is the characteristic of the residue
field of $A$, and $p\notin\fm_A^2$, we let $\bOmega_A$ be
the $k_A^{1/p}$-module introduced in \eqref{subsec_bOmega},
and if $p\in\fm_A^2$, we set
$$
\bOmega_A:=\Omega^1_{A/\Z}\otimes_Ak_A^{1/p}.
$$
To any local ring homomorphism $\phi:A\to B$ we attach
the $k_A^{1/p}$-linear map
$$
\bOmega_\phi:\bOmega_A\to\bOmega_B
$$
defined as follows :

\begin{itemize}
\item
In case $p\notin\fm_B^2$, then obviously $p\notin\fm_A^2$.
Denote by $\bar\phi:k_A^{1/p}\to k_B^{1/p}$ the map
induced by $\phi$. By \eqref{subsec_force-nextpage},
we have an induced $A$-linear map
$\bar\Omega_\phi\otimes_A\bar\phi:
\bar\Omega_A\otimes_Ak_A^{1/p}\to\bar\Omega_B\otimes_Bk_B^{1/p}$.
Also, the naturality of the ghost component map
$\bar\omega_0$ yields a $k_A$-linear map
$V_1(\phi):V_1(k_A)\to V_1(k_B)$, and indeed, $V_1(\phi)$
is even $k_A^{1/p}$-linear. Moreover, for every $a\in k_A$
and $b\in k_A^{1/p}$ we may compute
$$
\begin{aligned}
\bpsi_B(\bar\Omega_\phi(j_A(0,a))\otimes\bar\phi(b))
=\, & \bpsi_B(j_B(V_1(\phi)(0,a))\otimes\bar\phi(b)) \\
=\, & \bj_B\circ\psi_B(V_1(\phi)(0,a)\otimes\bar\phi(b)) \\
=\, & \bj_B\circ V_1(\phi)\circ\psi_A((0,a)\otimes b)
\end{aligned}
$$
so the pair $(\bpsi_B\circ(\bar\Omega_\phi\otimes\bar\phi),
\bj_B\circ V_1(\phi))$ determines a unique $k_A^{1/p}$-linear
map $\bOmega_\phi:\bOmega_A\to\bOmega_B$ as required.
\item
In case $p\in\fm^2_A$, then $p\in\fm^2_B$, and
we set $\bOmega_\phi:=\Omega^1_\phi\otimes_A\bar\phi$,
where $\Omega^1_\phi:\Omega^1_{A/\Z}\to\Omega^1_{B/\Z}$ is
the natural $A$-linear map.
\item
Lastly, if $p\notin\fm_A^2$ but $p\in\fm_B^2$, we
let $\bOmega_\phi$ be the composition of the natural
projection $\bOmega_A\to\Omega^1_A\otimes_Ak_A^{1/p}$
and $\Omega^1_\phi\otimes_A\bar\phi$.
\end{itemize}
If $\phi':B\to C$ is another morphism in $\mathbf{Local}$,
it is easy to verify that
$\bOmega_{\phi'\circ\phi}=\bOmega_{\phi'}\circ\bOmega_\phi$,
so we do get a functor as sought (details left to the
reader).

\sset\subsubsection{}\label{subsec_same-vein}
In the same vein, we may extend the map $\bd_A$
to a natural transformation of functors. Namely, if
$p\notin\fm_A^2$, we let $\bd_A$ be the map
introduced in \eqref{subsec_bOmega}, and if $p\in\fm_A^2$,
we let $\bd_A:A\to\bOmega_A$ be the map induced by
the universal derivation $A\to\Omega^1_{A/\Z}$. Then
it is easily seen that we have a commutative diagram
$$
\xymatrix{
A \ar[r]^-{\bd_A} \ar[d]_\phi &
\bOmega_A \ar[d]^{\bOmega_\phi} \\
B \ar[r]^-{\bd_B} & \bOmega_B
}$$
for every local ring homomorphism $\phi:A\to B$
(details left to the reader).

\begin{proposition}\label{prop_viens-viens}
Let $(A,\fm_A)$, $(B,\fm_B)$ be any two local rings,
$\phi:A\to B$ a local ring homomorphism, and denote
by $p$ the residue characteristic of $A$ and $B$.
We have :
\begin{enumerate}
\item
Suppose that $\phi$ is formally smooth for the $\fm_A$-adic
topology on $A$ and the $\fm_B$-adic topology on $B$. Then
$\phi$ and $\bOmega_\phi$ induce injective maps
$$
(\fm_A/\fm_A^2)\otimes_{k_A}k_B\to\fm_B/\fm_B^2
\qquad
\bOmega_A\otimes_{k^{1/p}_A}k^{1/p}_B\to\bOmega_B.
$$
\item
Suppose that
\begin{enumerate}
\item
$\fm_A\cdot B=\fm_B$.
\item
The induced residue field extension $k_A\to k_B$ is algebraic
and separable.
\item
$\phi$ is flat.
\end{enumerate}
Then $\bOmega_\phi$ induces an isomorphism of
$k_B^{1/p}$-vector spaces :
$$
\bOmega_A\otimes_AB\isom\bOmega_B.
$$
\item
If $B=A/\fm_A^2$ and $\phi$ is the natural
projection, then $\bOmega_\phi$ is an isomorphism.
\item
The functor $\bOmega_\bullet$ and the natural transformation
$\bd_\bullet$ commute with all filtered colimits.
\end{enumerate}
\end{proposition}
\begin{proof}(ii): Assumption (a) and the flatness of
$\phi$ imply that the induced map
$$
(\fm_A/\fm_A^2)\otimes_{k_A}k_B\to\fm_B/\fm_B^2
$$
is an isomorphism (\cite[Th.22.3]{Mat}), so $p\in\fm_A^2$
if and only if $p\in\fm_B^2$. On the other hand, under
assumption (b), we have
\set\begin{equation}\label{eq_same-perfect}
k^{1/p}_B=k_A^{1/p}\cdot k_B.
\end{equation}
Now, consider first the case where $p\in\fm_A^2$;
then $\bOmega_\phi=\Omega^1_\phi\otimes_A\bar\phi$
(notation of \eqref{subsec_bOmega-funct}), so
the assertion follows from \eqref{eq_same-perfect}
together with the following

\begin{claim}\label{cl_appealing}
If $p\in\fm^2_A$, then $\Omega^1_\phi$ induces a
$k_B$-linear isomorphism
$$
\Omega^1_{A/\Z}\otimes_Ak_B\to\Omega^1_{B/\Z}\otimes_Bk_B.
$$
\end{claim}
\begin{pfclaim} Under the standing assumptions,
\eqref{subsec_please} gives a commutative ladder with
short exact rows
\set\begin{equation}\label{eq_ladder-Kaehler}
{\diagram
0 \ar[r] & (\fm_A/\fm_A^2)\otimes_{k_A}k_B \ar[r] \ar[d] &
\Omega^1_{A/\Z}\otimes_Ak_B \ar[r] \ar[d] &
\Omega^1_{k_A/\Z}\otimes_{k_A}k_B \ar[r] \ar[d] & 0 \\
0 \ar[r] & \fm_B/\fm_B^2 \ar[r] &
\Omega^1_{B/\Z}\otimes_Bk_B \ar[r] &
\Omega^1_{k_B/\Z} \ar[r] & 0
\enddiagram}
\end{equation}
and we have already remarked that the left vertical
arrow is an isomorphism; the same holds for the right
vertical arrow, due to assumption (b). The claim
follows.
\end{pfclaim}

The case where $p\notin\fm^2_B$ is treated likewise :
one argue as in the proof of claim \ref{cl_appealing},
except that, instead of appealing to \eqref{eq_finally-exact},
one invokes proposition \ref{prop_lift-lemma} and the
naturality of $\bar\bd_\bullet$ established in
\eqref{subsec_same-vein} : the details shall be
left to the reader.

(iii): For the case where $p\in\fm_A^2$, we consider again
the induced ladder \eqref{eq_ladder-Kaehler} : under the
standing assumptions, clearly the first and third vertical
arrows are isomorphisms, so the same holds for the middle one.

For the case where $p\notin\fm_A^2$, one argues again
likewise : the ladder \eqref{eq_ladder-Kaehler} is
replaced by the corresponding ladder for $\bOmega_\bullet$,
provided by proposition \ref{prop_lift-lemma} and the
naturality of $\bar\bd_\bullet$ : the details shall
be left to the reader.

(iv): Let $((A_\lambda,\fm_\lambda)~|~\lambda\in\Lambda)$ be a
filtered system of local rings and local ring homomorphisms.
Notice that if $p\in\fm_\lambda^2$ for some
$\lambda\in\Lambda$, then $p\in\fm_\mu^2$ for every
$\mu\geq\lambda$ in $\Lambda$, and since $\Lambda$
is filtered, we may then replace $\Lambda$ with
$\Lambda/\lambda$, and assume that $p\in\fm_\mu$ for
every $\mu\in\Lambda$. In this case, the contention
follows straightforwardly from the corresponding
assertion for the functor $\Omega^1_\bullet$ and the
universal derivation $d_\bullet$.

It remains to consider the case where $p\notin\fm_\lambda$
for every $\lambda\in\Lambda$. For this, one applies
proposition \ref{prop_lift-lemma}, which reduces
again the contention to the previous case (details
left to the reader).

(i): Let us begin with the following general result :

\begin{claim}\label{cl_Cartier}
Let $K$ be any field, and $L$ a field extension of $K$ of
finite type. Then
$$
\dim_L\Omega^1_{L/K}-\dim_LH_1\L_{L/K}=\tr.\deg[L:K].
$$ 
\end{claim}
\begin{pfclaim} This is known as {\em Cartier's identity},
and a proof is given in \cite[Ch.0, Th.21.7.1]{EGAIV}. We
present a proof via the cotangent complex formalism.
Suppose first that $K$ is of finite type over its prime
field $K_0$ (so $K_0$ is either $\Q$ or a finite field of
prime order). In this case, notice that $H_1\L_{L/K_0}=0$
(cp. the proof of lemma \ref{lem_finally-exact}); then
the transitivity triangle for the sequence of maps
$K_0\to K\to L$ (\cite[Ch.II, \S2.1.2.1]{Il}) yields an
exact sequence
$$
0\to H_1\L_{L/K}\to\Omega^1_{K/K_0}\otimes_KL\to\Omega^1_{L/K_0}
\to\Omega^1_{L/K}\to 0
$$
and the claim follows easily, after one remarks that
$\dim_K\Omega^1_{K/K_0}=\tr.\deg[K:K_0]$, and likewise
for $\Omega^1_{L/K_0}$. Next, let $K$ be an arbitrary
field, and write $K$ as the union of the filtered family
$(K_\lambda~|~\lambda\in\Lambda)$ of its subfields that
are finitely generated over $K_0$. Choose elements
$x_1,\dots,x_t\in L$ such that $K(x_1,\dots,x_t)=L$. Let
$\phi:K[X_1,\dots,X_t]\to L$ be the map of $K$-algebras
given by the rule : $X_i\mapsto x_i$ for $i=1,\dots,t$.
Then $I:=\Ker\,\phi$ is a finitely generated ideal, so we
may find $\lambda\in\Lambda$ and an ideal
$I_\lambda\subset K_\lambda[X_1,\dots,X_t]$ such that
$I=I_\lambda\otimes_{K_\lambda}K$ (details left to the
reader). Denote by $L_\lambda$ the field of fractions of
$K_\lambda[T_1,\dots,T_t]/I_\lambda$; it is easily seen
that $L_\lambda\otimes_{K_\lambda}K$ is an integral
domain, and its field of fractions is a $K$-algebra
naturally isomorphic to $L$. Especially, we have
$\tr.\deg[L_\lambda:K_\lambda]=\tr.\deg[L:K]$, and on
the other hand, there is a natural isomorphism in
$\sD^-(L\Mod)$ :
$$
\L_{L_\lambda/K_\lambda}\otimes_{L_\lambda}L\isom\L_{L/K}
$$
(\cite[Ch.II, Prop.2.2.1, Cor.2.3.1.1]{Il}). Then the sought
identity for the extension $K\subset L$ follows from the same
identity for the extension $K_\lambda\subset L_\lambda$. The
latter is already known, by the previous case.
\end{pfclaim}

Next, we consider the following special case :

\begin{claim}\label{cl_yet-unnamed}
Let $n\in\N$ be any integer, and $\fq\subset A[T_1,\dots,T_n]$
any prime ideal such that $\fq\cap A=\fm_A$. Set
$R:=A[T_1,\dots,T_n]_\fq$ and denote by $\fm_R$ and $k_R$
respectively the maximal ideal and the residue field of
$R$. Then :
\begin{enumerate}
\item
$\dim_{k_R}\fm_R/\fm_R^2=n+\dim_{k_A}\fm_A/\fm^2_A$.
\item
The natural map
$\gamma:(\fm_A/\fm^2_A)\otimes_{k_A}k_R\to\fm_R/\fm_R^2$
is injective.
\item
The natural map $\Omega^1_{A/\Z}\otimes_AR\to\Omega^1_{R/\Z}$
is injective, and its image is a direct summand of
$\Omega^1_{R/\Z}$.
\item
The natural map $\bOmega_A\otimes_{k_A^{1/p}}k_R^{1/p}\to\bOmega_R$
is injective.
\end{enumerate}
\end{claim}
\begin{pfclaim}(i): According to \cite[Th.6.5.12(i)]{Ga-Ra}, we
have $H_i\L_{k_R/k_A}=0$ for every $i>1$ (we apply {\em loc.cit.}
to the extension $(k_A,|\cdot|_{k_A})\subset(k_R,|\cdot|_{k_B})$
of valued fields with trivial valuations); then the transitivity
triangle for the cotangent complex relative to the maps
$A\to k_A\to k_R$ yields a short exact sequence of $k_R$-vector
spaces :
$$
0\to H_1\L_{k_A/A}\otimes_{k_A}k_R\xrightarrow{\ \alpha\ }
H_1\L_{k_R/A}\to H_1\L_{k_R/k_A}\to 0
$$
(\cite[Ch.II, \S2.1.2.1]{Il}). Likewise, since
$H_i\L_{R/A}=0$ for every $i>1$ (\cite[Ch.II, Cor.1.2.6.3]{Il}),
the sequence of maps $A\to R\to k_R$ yields an exact sequence
of $k_R$-vector spaces :
$$
0\to H_1\L_{k_R/A}\xrightarrow{\ \beta\ }
H_1\L_{k_R/R}\to\Omega^1_{R/A}\otimes_Rk_R\to\Omega_{k_R/A}\to 0.
$$
However, we have natural isomorphisms
$$
H_1\L_{k_A/A}\isom\fm_A/\fm_A^2
\qquad
H_1\L_{k_R/R}\isom\fm_R/\fm_R^2
$$
(\cite[Ch.III, Cor.1.2.8.1]{Il}), and clearly
$\dim_{k_R}\Omega_{R/A}\otimes_Rk_R=n$. Thus :
$$
\begin{aligned}
\dim_{k_A}\fm_A/\fm_A^2=\: &
\dim_{k_R}H_1\L_{k_R/A}-\dim_{k_R}H_1\L_{k_R/k_A} \\
\dim_{k_R}\fm_R/\fm_R^2=\: &\dim_{k_R}H_1\L_{k_R/A}+n-
\dim_{k_R}\Omega_{k_R/k_A}.
\end{aligned}
$$
Taking into account the identity
$$
\dim_{k_R}H_1\L_{k_R/k_A}=\dim_{k_R}\Omega_{k_R/k_A}
$$
provided by claim \ref{cl_Cartier}, the assertion follows.

(ii): Notice that the composition of $\alpha$ and $\beta$
yields an injective map
$(\fm_A/\fm^2_A)\otimes_{k_A}k_R\to\fm_R/\fm_R^2$ so it suffices
to check that this composition equals $\gamma$. However, let
$$
\Sigma
\quad :\quad
H_1\L_{k_R/R}\xrightarrow{\ d\ } H_0k_R\otimes_R\L_{R/A}
$$
be the boundary map of the transitivity triangle
\set\begin{equation}\label{eq_first-transit}
k_R\otimes_R\L_{R/A}\to\L_{k_R/A}\to\L_{k_R/R}\to
k_R\otimes_R\sigma\L_{R/A}
\end{equation}
arising from the sequence $A\to R\to k_R$; we regard $\Sigma$
as a complex placed in degrees $[-1,0]$, and then it is clear
that the triangle $\tau^{\geq -1}\eqref{eq_first-transit}$ is
naturally isomorphic in $\sD(k_R\Mod)$ to the triangle
\set\begin{equation}\label{eq_oh-no-another}
k_R\otimes_RH_0\L_{R/A}[0]\to\Sigma\to H_1\L_{k_R/R}[1]\to
k_R\otimes_RH_0\L_{R/A}[1].
\end{equation}
According to \cite[III.1.2.9.1]{Il}, the complex $\Sigma$
is naturally isomorphic to the complex
$$
\Theta
\quad:\quad
\fm_R/\fm_R^2\xrightarrow{\ -d_{R/A}\ }k_R\otimes_R\Omega^1_{R/A}
$$
(placed in degrees $[-1,0]$) where $d_{R/A}$ is induced
by the universal derivation $R\to\Omega^1_{R/A}$. Likewise,
we have natural isomorphisms
\set\begin{equation}\label{eq_identify}
\tau^{\geq -1}\L_{k_R/R}\isom\fm_R/\fm_R^2[1]
\qquad
\tau^{\geq -1}\L_{R/A}\isom\Omega^1_{R/A}[0]
\end{equation}
and under these identifications, \ref{eq_oh-no-another} is
the obvious triangle deduced from $\Theta$. Especially,
the map $\beta$ is naturally identified with the inclusion
map $\Ker\,d_{R/A}\to\fm_R/\fm_R^2$.

Likewise, there exists a natural identification
\set\begin{equation}\label{eq_injective-two}
\tau^{\geq -1}\L_{k_A/A}\isom\fm_A/\fm^2_A[1]
\qquad
\text{in $\sD(k_A\Mod)$}
\end{equation}
as well as a natural map of complexes
\set\begin{equation}\label{eq_injective-one}
k_R\otimes_{k_A}(\fm_A/\fm^2_A)[1]\to\Theta
\end{equation}
deduced from $\gamma$. To conclude, it suffices to check
that the map
\set\begin{equation}\label{eq_tau_-and-alpha}
\tau^{\geq -1}k_R\otimes_A\L_{k_A/A}\to\tau^{\geq -1}\L_{k_R/A}
\end{equation}
coming from the transitivity triangle for the sequence
$A\to k_A\to k_R$, corresponds to the morphism
\eqref{eq_injective-one}, under the
identification \eqref{eq_injective-two} and the previous
identification of $\tau^{\geq -1}\L_{k_R/A}$ with $\Theta$.
To this aim, it suffices to compare the maps obtained by
applying to these two morphisms the functor $\Ext^1_{k_R}(-,M[0])$,
for arbitrary $k_R$-modules $M$. Now, recall that there
exists a natural isomorphism
$$
\Ext^1_{k_R}(\tau^{\geq -1}\L_{k_R/R},M[0])\isom\Exal_R(k_R,M)
\qquad
\text{for every $k_R$-module $M$}.
$$
As explained in \cite[III.1.2.8]{Il}, under the identification
\eqref{eq_identify}, this becomes the following $k_R$-linear
isomorphism
\set\begin{equation}\label{eq_realize}
\Hom_{k_R}(\fm_R/\fm_R^2,M)\isom\Exal_R(k_R,M)
\qquad
\phi\mapsto \phi*U
\end{equation}
(notation of \cite[\S2.5.5]{Ga-Ra}), where
$$
U
\quad : \quad
0\to\fm_R/\fm_R^2\to R/\fm_R^2\to k_R\to 0
$$
is the natural extension. There is a natural surjection
$$
\Hom_{k_R}(\fm_R/\fm_R^2,M)\to\Ext^1_{k_R}(\Theta,M[0])
$$
and the foregoing implies that the natural isomorphism
$$
\Ext^1_{k_R}(\tau^{\geq -1}\L_{k_R/A},M[0])\isom\Exal_A(k_R,M)
$$
is also realized as in \eqref{eq_realize} : given a class $c$
in $\Ext^1_{k_R}(\Theta,M|0])$, take an arbitrary representative
$\phi:\fm_R/\fm_R^2\to M$, and the correspondence associates with
$c$ the class of the push out $\phi*U$.

By the same token, we have a natural isomorphism
$$
\Ext^1_{k_R}(\tau^{\geq -1}k_R\otimes_{k_A}\L_{k_A/A},M(0])
\isom\Exal_A(k_A,M)
\qquad
\text{for every $k_R$-module $M$}
$$
and on the one hand, the map
$\Ext^1_{k_R}(\eqref{eq_tau_-and-alpha},M[0])$ is identified
naturally with the map
$$
\Exal_A(k_R,M)\to\Exal_A(k_A,M)
$$
given by pull back along the inclusion map $\iota:k_A\to k_R$.
On the other hand, by \cite[III.1.2.8]{Il}, the identification
 \eqref{eq_identify} induces the isomorphism
$$
\Hom_{k_A}(\fm_A/\fm_A^2,M)\isom\Exal_A(k_A,M)
\qquad
\phi'\mapsto\phi'*U'
$$
where
$$
U'
\quad :\quad
0\to\fm_A/\fm^2_A\to A/\fm^2_A\to k_A\to 0
$$
is the natural extension. So finally, the assertion boils
down to the identity
$$
(\phi\circ\gamma)*U'=\phi*U*\iota
\qquad
\text{in $\Exal_A(k_A,M)$}
$$
for every $k_R$-module $M$ and every $\phi:\fm_R/\fm_R^2\to M$.
We leave the verification as an exercise for the reader. 

(iii) is a standard calculation (more
precisely, the complement of $\Omega^1_{A/\Z}\otimes_AR$
in $\Omega^1_{R/\Z}$ is the free $R$-module generated by
$dT_1,\dots,dT_n$).

(iv): If $p\in\fm^2_A$, the assertion follows from (iii).
Suppose then that $p\notin\fm^2_A$. By inspecting the
constructions, we get a natural commutative ladder of
$k_R^{1/p}$-vector spaces
$$
\xymatrix{
0 \ar[r] & V_1(k_A)\otimes_{k_A^{1/p}}k_R^{1/p} \ar[r] \ar[d] &
\bOmega_A\otimes_{k_A^{1/p}}k_R^{1/p} \ar[r] \ar[d] &
\Omega^1_{A/\Z}\otimes_Ak_R^{1/p} \ar[r] \ar[d] & 0 \\
0 \ar[r] & V_1(k_R) \ar[r] & \bOmega_R \ar[r] &
\Omega^1_{R/\Z}\otimes_Rk_R^{1/p} \ar[r] & 0
}$$
whose central vertical arrow is the map of the claim.
However, it is easily seen that the left vertical arrow
is an isomorphism, so the assertion follows from (iii).
\end{pfclaim}

Now, pick a free polynomial $A$-algebra $C$
and a surjective map $C\to B$ of $A$-algebras; let
$\fq\subset C$ be the preimage of $\fm_B$, set
$R:=C_\fq$, and let $\fm_R$ be the maximal ideal
of $R$. Since $B$ is formally smooth over $A$ for
the topologies of the maximal ideals, the induced
surjection $R/\fm_R^2\to B/\fm_B^2$ admits a
section $B/\fm_B^2\to R/\fm_R^2$ which is also a
local map of $A$-algebras. In view of (iii), there
follow commutative diagrams of $k_B$-linear and
respectively $k_B^{1/p}$-linear maps
$$
\xymatrix{ (\fm_A/\fm_A^2)\otimes_{k_A}k_B
\ar[rd] \ar[rr] & & \fm_B/\fm_B^2 \ar[ld] &
\bOmega_A\otimes_{k_A}k_B^{1/p} \ar[rr] \ar[rd]
& & \bOmega_B \ar[ld] \\
& \fm_R/\fm_R^2 & & & \bOmega_R
}$$
which reduce the assertion to the case where $B=R$.
In this case, write $C$ as the filtered colimit of
the system $(C_\lambda~|~\lambda\in\Lambda)$ of its
free polynomial $A$-subalgebras of finite type;
for each $\lambda\in\Lambda$, let
$\fq_\lambda:=\fq\cap C_\lambda$ and
$R_\lambda:=(C_\lambda)_{\fq_\lambda}$. Clearly $R$ is
the filtered colimit of the system
$(R_\lambda~|~\lambda\in\Lambda)$, so the assertion
follows from (iv) and claim \ref{cl_yet-unnamed}(ii).
\end{proof}

\sset\subsubsection{}\label{subsec_setup-criterion}
Keep the notation of \eqref{subsec_please}, and let now
$f_1,\dots,f_n$ be a finite sequence of elements of $A$,
and $e_1,\dots,e_n$ a sequence of integers such that
$e_i>1$ for every $i=1,\dots,n$. Set
$$
C:=A[T_1,\dots,T_n]/(T_1^{e_1}-f_1,\dots,T^{e_n}_n-f_n).
$$
Fix a prime ideal $\fn\subset C$ such that $\fn\cap A=\fm_A$,
and let $B:=C_\fn$. So the induced map $A\to B$ is a local
ring homomorphism; we denote by $\fm_B$ the maximal ideal
of $B$, and set $k_B:=B/\fm_B$. Also, let $\nu:=\dim_{k_A}E$,
where $E\subset\bOmega_A$ is the $k^{1/p}_A$-vector space
spanned by $\bd_Af_1,\dots,\bd_Af_n$.

\begin{theorem}\label{th_regu-criterion}
In the situation of \eqref{subsec_setup-criterion}, suppose
moreover that :
\begin{enumerate}
\alphaenu
\item
$f_i\in\fm_A$, for every $i\leq n$ such that $p$ does not divide
$e_i$.
\item
$\fm_A/\fm_A^2$ is a finite dimensional $k_A$-vector space.
\end{enumerate}
Then $\fm_B/\fm_B^2$ is a finite dimensional $k_B$-vector
space, and we have :
$$
\dim_{k_B}\fm_B/\fm_B^2=n+\dim_{k_A}\fm_A/\fm_A^2-\nu.
$$
\end{theorem}
\begin{proof} Let $\fq\subset A[T_1,\dots,T_n]$ be the
preimage of $\fn$, set $R:=A[T_1,\dots,T_n]_\fq$, and
denote by $\fp$ the maximal ideal of $R$, and by
$F\subset\fp/\fp^2$ the $k_B$-vector space spanned
by $T_1^{e_1}-f_1,\dots,T_n^{e_n}-f_n$. Clearly,
we have a short exact sequence
$$
0\to F\to\fp/\fp^2\to\fm_B/\fm_B^2\to 0.
$$
Suppose first that $p\notin\fm_A^2$, in which case
$p\notin\fp^2$, by claim \ref{cl_yet-unnamed}(ii), hence
$\bOmega_R$ is defined by \eqref{subsec_bOmega}. On the
other hand, notice that
$\bd(T_i^{e_i})=\bd(T_i^{e_i}-f_i)+\bd(f_i)$ in $\bOmega_R$,
since $T_i^{e_i}-f_i\in\fp$. Now, if $e_i$ is a multiple of
$p$, Leibniz's rule yields
$\bd(T_i^{e_i})=e_i\cdot T_i^{e_i-1}\bd(T_i)=0$.
If $e_i$ is not a multiple of $p$, we have $f_i\in\fm_A$
by assumption; hence $T_i\in\fp$ and therefore $\bd(T_i^{e_i})=0$
again, since $e_i>1$. In either case, we find
$$
\bd(f_i)=-\bar\bd(T_i^{e_i}-f_i)
\qquad
\text{in $\bOmega_R$, for every $i=1,\dots,n$.}
$$
In view of proposition \ref{prop_lift-lemma}, it follows that
$$
\dim_{k_B}\fm_B/\fm^2_B=\dim_{k_B}\fp/\fp^2-\dim_{k^{1/p}_B}E'
$$
where $E'\subset\bOmega_R$ is the $k^{1/p}_B$-vector space
spanned by $\bd(f_1),\dots,\bd(f_n)$. Then the assertion
follows from claim \ref{cl_yet-unnamed}(i,iv).

Lastly, suppose $p\in\fm_A^2$, so that $p\in\fp^2$ as well.
Arguing as in the foregoing case, we see that $\dim_{k_B}F$
equals the dimension of the $k_B$-vector subspace of
$\Omega^1_{R/\Z}\otimes_Rk_B$ spanned by $df_1,\dots,df_n$,
and in view of claim \ref{cl_yet-unnamed}(iii), the
latter equals $\nu$, whence the contention. 
\end{proof}

\begin{corollary}\label{cor_regu-criterion}
In the situation of theorem {\em\ref{th_regu-criterion}}, the
following conditions are equivalent :
\begin{enumerate}
\alphaenu
\item
$A$ is a regular local ring, and $\nu=n$.
\item
$B$ is a regular local ring.
\end{enumerate}
\end{corollary}
\begin{proof} Suppose first that (b) holds. Since the
map $A\to B$ is faithfully flat, it is easily seen that
$A$ is noetherian, and then \cite[Ch.0, Prop.17.3.3(i)]{EGAIV}
shows already that $A$ is a regular local ring. Moreover,
$C$ is clearly a finite $A$-algebra, therefore
$\dim_{k_B}\fm_B/\fm_B^2=\dim B\leq\dim A=\dim_{k_A}\fm_A/\fm_A^2$.
Theorem \ref{th_regu-criterion} then implies that $\nu=n$, so
(a) holds.

Next, suppose that (a) holds. We apply theorem
\ref{th_regu-criterion} as in the foregoing, to deduce that
$\dim B=\dim A=\dim_{k_B}\fm_B/\fm_B^2$, whence (b) : details
left to the reader.
\end{proof}

\begin{remark}\label{rem_afterthought}
(i)\ \
Keep the notation of corollary \ref{cor_regu-criterion}.
In \cite[Ch.0, Th.22.5.4]{EGAIV} it is asserted that condition
(b) is equivalent to the following :

(a')\ \ $A$ is regular and the space
$E\subset\Omega^1_{A/\Z}\otimes_Ak_A$ spanned by $df_1,\dots,df_n$
has dimension $n$.

When $p\in\fm_A^2$, this condition (a') agrees with our condition
(a), so in this case of course we do have (a')$\Leftrightarrow$(b).
However, the latter equivalence fails in general, in case
$p\notin\fm_A$ : the mistake is found in
\cite[Ch.0, Rem.22.4.8]{EGAIV}, which is false. The implication
(a')$\Rightarrow$(b) does remain true in all cases : this is
easily deduced from theorem \ref{th_regu-criterion}, since
the image of $\bd_A(f)$ in $\Omega^1_{A/\Z}\otimes_{k_A}k_A^{1/p}$
agrees with $df\otimes 1$, for every $f\in A$. The proof of
{\em loc.cit.} is correct for $p\in\fm_A^2$, which is the only
case that is used in the proof of corollary \ref{cor_regu-criterion}.

(ii)\ \
Moreover, in the situation of corollary \ref{cor_regu-criterion},
suppose that $B$ is a regular local ring. Then the sequence
$$
(\bd_Bf^{1/e_1}_1,\dots,\bd_Bf^{1/e_n}_n)
$$
spans a $k^{1/p}_B$-vector subspace of $\bOmega_B$ of dimension $n$.
Indeed, consider the ring
$$
C:=B[T_1,\dots,T_n]/(T_1^p-f^{1/e_1}_1,\dots,T^p_n-f^{1/e_n}_n).
$$
Corollary \ref{cor_regu-criterion} applies to the extension
$A\subset C$, the sequence $(f_1,\dots,f_n)$, and the
sequence of integers $(pe_1,\dots,pe_n)$, so $C$ is a
regular local ring. But the same corollary applies as well
to the extension $B\subset C$, the sequence
$(f^{1/e_1}_1,\dots,f_n^{1/e_n})$, and the sequence of integers
$(p,\dots,p)$, and yields the assertion.

(iii)\ \
Let us say that the sequence $(f_1,\dots,f_n)$ is
{\em maximal in $A$} if
$$
(\bd_Af_1,\dots,\bd_Af_n)
$$
is a basis of the $k_A^{1/p}$-vector space $\bOmega_A$.
Then, in the situation of (ii), we claim that the sequence
$(f_1,\dots,f_n)$ is maximal in $A$ if and only if the
sequence $(f_1^{1/e_1},\dots,f_n^{1/e_n})$ is maximal in $B$.
Indeed, under the current assumptions we have
$\dim_{k_A}\Omega^1_{k_A/\Z}=\dim_{k_B}\Omega^1_{k_B/\Z}$, since
these integers are equal to the transcendence degree of $k_A$
(and $k_B$) over $\F_p$, and on the other hand
$\dim_{k_A}\fm_A/\fm_A^2=\dim_{k_B}\fm_B/\fm_B^2$, since $A$ and
$B$ are regular local rings of the same dimension; then the
assertion follows from (ii), lemma \ref{lem_finally-exact} and
proposition \ref{prop_lift-lemma} (details left to the reader).

(iv)\ \
Likewise, consider any morphism $\phi:A\to B$ in $\mathbf{Local}$
fulfilling conditions (a),(b) and (c) of proposition
\ref{prop_viens-viens}(ii). Then clearly the sequence
$(f_1,\dots,f_n)$ is maximal in $A$ if and only if the sequence
$(\phi(f_1),\dots,\phi(f_n))$ is maximal in $B$.
\end{remark}

\sset\subsubsection{}\label{subsec_again-Frobby}
Let $p>0$ be a prime integer, and $A$ an $\F_p$-algebra.
Denote by $\Phi_A:A\to A$ the {\em Frobenius endomorphism\/}
of $A$, given by the rule : $a\mapsto a^p$ for every $a\in A$.
For every $A$-module $M$, we let $M_{(\Phi)}$ be the $A$-module
obtained from $M$ via restriction of scalars along the map
$\Phi_A$ (that is, $a\cdot m:=a^pm$ for every $a\in A$ and
$m\in M$). Notice that $\Phi_A$ is an $A$-linear map
$A\to A_{(\Phi)}$. Theorem \ref{th_Kunz-by-Matsu}, and part (i)
of the following theorem \ref{th_Kunz-exc} are due to E.Kunz.

\begin{theorem}\label{th_Kunz-by-Matsu}
Let $A$ be a noetherian local $\F_p$-algebra.
Then the following conditions are equivalent :
\begin{enumerate}
\item
$A$ is regular.
\item
$\Phi_{\!A}$ is a flat ring homomorphism.
\item
There exists $n>0$ such that $\Phi_{\!A}^n$ is a flat ring
homomorphism.
\end{enumerate}
\end{theorem}
\begin{proof}(i)$\Rightarrow$(ii): Let $A^\wedge$ be the
completion of $A$, and $f:A\to A^\wedge$ the natural map.
Clearly
$$
f\circ\Phi_A=\Phi_{A^\wedge}\circ f.
$$
Since $f$ is faithfully flat, it follows that $\Phi_A$ is
flat if and only if the same holds for $\Phi_{A^\wedge}$,
so we may replace $A$ by $A^\wedge$, and assume from start
that $A$ is complete, hence $A=k[[T_1,\dots,T_d]]$, for
a field $k$ of characteristic $p$, and $d=\dim A$
(\cite[Ch.0, Th.19.6.4]{EGAIV}). Then, it is easily seen that
$$
\Phi_A(A)=A^p=k^p[[T_1^p,\dots,T^p_d]].
$$
Set $B:=k[[T_1^p,\dots,T^p_d]]$; the ring $A$ is a free
$B$-module (of rank $p^d$), hence it suffices to check that
the inclusion map $A^p\to B$ is flat. However, denote
by $\fm$ the maximal ideal of $A^p$; clearly $B$ is an
$\fm$-adically ideal-separated $A$-module (see
\cite[p.174, Def.]{Mat}), hence it suffices to check that
$B/\fm^k B$ is a flat $A^p/\fm^k$-module for every $k>0$
(\cite[Th.22.3]{Mat}). The latter is clear, since
$k[T_1^p,\dots,T_d^p]$ is a flat $k^p[T_1^p,\dots,T^p_d]$-module.

(ii)$\Rightarrow$(iii) is obvious.

(iii)$\Rightarrow$(i): Notice first that $\Spec\,\Phi^n_A$
is the identity map on the topological space underlying $\Spec\,A$;
especially, $\Phi^n_A$ is flat if and only if it is faithfully
flat, and the latter condition implies that $\Phi^n_A$ is
injective. We easily deduce that if (iii) holds, then $A$
is reduced. Now, consider quite generally, any finite
system $x_\bullet:=(x_1,\dots,x_t)$ of elements of $A$, and
let $I\subset A$ be the ideal generated by $x_\bullet$; we
shall say that $x_\bullet$ is a system of {\em independent\/}
elements, if $I/I^2$ is a free $A/I$-module of rank $n$.
We show first the following :

\begin{claim}\label{cl_Leech}
Let $y,z,x_2,\dots,x_t$ be a family of elements of $A$, such
that $x_\bullet:=(yz,x_2,\dots,x_t)$ is a system of independent
elements, and denote by $J\subset A$ the ideal generated by
$x_\bullet$. We have :
\begin{enumerate}
\item
$(y,x_2,\dots,x_t)$ is a system of independent elements of $A$.
\item
If $\length_AA/J$ is finite, then
$$
\length_AA/J=
\length_AA/(y,x_2,\dots,x_t)+\length_AA/(z,x_2,\dots,x_t).
$$
\end{enumerate}
\end{claim}
\begin{pfclaim}(i): Suppose $a_1y+a_2x_2+\cdots+a_tx_t=0$ is a
linear relation with $a_1,\dots,a_t\in A$, and let $I\subset A$
be the ideal generated by $y,x_2,\dots,x_t$. We have to show that
$a_1,\dots,a_t\in I$. However, as $a_1yz+a_2zx_2+\cdots+a_tzx_t=0$,
it follows by assumption, that $a_1$ lies in $J\subset I$. Write
$a_1=b_1yz+b_2x_2+\cdots+b_tx_t$; then
$$
b_1y^2z+(a_2+b_2y)x_2+\cdots+(a_t+b_ty)x_t=0.
$$
Therefore $a_i+b_iy\in J$ for $i=2,\dots,t$, and therefore
$a_2,\dots,a_t\in I$, as required.

(ii): It suffices to show that the natural map of $A$-modules :
$$
A/(z,x_2,\cdots,x_t)\to I/J
\qquad
a\mapsto ay+J
$$
is an isomorphism. However, the surjectivity is immediate. To
show the injectivity, suppose that $ay\in J$, {\em i.e.}
$ay=b_1yz+b_2x_2+\cdots+b_tx_t$ for some $b_1,\dots,b_t\in A$;
we deduce that $(b_1z-a)yz+b_2zx_2+\cdots+b_tzx_t=0$, hence
$a-b_1z\in J$ by assumption, so $a$ lies in the ideal generated
by $z,x_2,\dots,x_t$, as required.
\end{pfclaim}

Now, set $q:=p^n$, and $A_\nu:=A^{q^\nu}\subset A$ for every
integer $\nu>0$; pick a minimal system $x_\bullet:=(x_1,\dots,x_t)$
of generators of the maximal ideal $\fm_A$ of $A$, and notice that,
since $A$ is reduced, $\Phi^{n\nu}_A$ induces an isomorphism
$A\to A_\nu$, hence $x_\bullet^{(\nu)}:=(x_1^{q^\nu},\dots,x_t^{q^\nu})$
is a minimal system of generators for the maximal ideal $\fm_\nu$
of $A_\nu$. Set as well $I_\nu:=\fm_\nu A$; since the inclusion
map $A_\nu\to A$ is flat by assumption for every $\nu>0$, we have
a natural isomorphism of $A$-modules
$$
(\fm_\nu/\fm_\nu^2)\otimes_{A_\nu}A\isom I_\nu/I_\nu^2.
$$
On the other hand, set $k_\nu:=A_\nu/\fm_\nu$; by Nakayama's lemma,
$\dim_{k_\nu}\fm_\nu/\fm_\nu^2=t$, so $I_\nu/I_\nu^2$ is a free
$A$-module of rank $t$, {\em i.e.} $x_\bullet^{(\nu)}$ is an
independent system of elements of $A$. From claim \ref{cl_Leech}(ii)
and a simple induction, we deduce that
\set\begin{equation}\label{eq_length-estim}
\length_AA/I_\nu=\length_{A^\wedge}A^\wedge/I_\nu A^\wedge=q^{\nu t}
\qquad
\text{for every $\nu>0$}
\end{equation}
(where the first equality holds, since $I_\nu$ is an open ideal
in the $\fm_A$-adic topology of $A$). According to
\cite[Ch.0, Th.19.9.8]{EGAIV} (and its proof), $A^\wedge$
contains a field isomorphic to $k_0:=A/\fm_A$, and the inclusion
map $k_0\to A^\wedge$ extends to a surjective ring homomorphism
$k_0[[X_1,\dots,X_t]]\to A^\wedge$, such that $X_i\mapsto x_i$
for $i=1,\dots,t$. Denote by $J$ the kernel of this surjection;
in view of \eqref{eq_length-estim}, we have
$$
\length_{A^\wedge}k_0[[X_1,\dots,X_t]]/(J,X_1^{q^\nu},\dots,X_t^{q^\nu})
=q^{\nu t}
$$
which means that $J\subset(X_1^{q^\nu},\dots,X_t^{q^\nu})$ for
every $\nu>0$. We conclude that $J=0$, and $A^\wedge=k_0[[X_1,\dots,X_t]]$
is regular, so the same holds for $A$.
\end{proof}

The last result of this section is a characterization of
regular local rings via the cotangent complex, borrowed from
\cite{An}, which shall be used in the following section on
excellent rings.

\begin{lemma}\label{lem_regu-seq-cot}
Let $A$ be a ring, $(a_1,\dots,a_n)$ a regular sequence
of elements of $A$, that generates an ideal $I\subset A$,
and set $A_0:=A/I$. Then there is a natural isomorphism
$$
\L_{A_0/A}\isom I/I^2[1]
\qquad
\text{in $\sD(A_0\Mod)$}
$$
and $I/I^2$ is a free $A_0$-module of rank $n$.
\end{lemma}
\begin{proof} Notice that $H_0\L_{A_0/A}=0$, and there is
a natural isomorphism $H_1\L_{A_0/A}\isom I/I^2$
(\cite[Ch.III, Cor.1.2.8.1]{Il}). There follows
a natural morphism $\L_{A_0/A}\to I/I^2[1]$, and we shall show
more precisely, that this morphism is an isomorphism.
We proceed by induction on $n$. Hence, suppose first that $n=1$,
set $a:=a_1$, and $B:=A[T]$, the free polynomial $A$-algebra in
one variable; define a map of $A$-algebras $B\to A$ by the rule
$T\mapsto a$. Set also $B_0:=B/TB$ (so $B_0$ is isomorphic
to $A$). Since $a$ is regular, it is easily seen that the
natural morphism
$$
B_0\derotimes_BA\to B_0\otimes_BA=A_0
$$
is an isomorphism in $\sD(B\Mod)$. It follows that the induced
morphism $\L_{B_0/B}\otimes_{B_0}A_0\to\L_{A_0/A}$ is an isomorphism
in $\sD(A_0\Mod)$ (\cite[Ch.II, Prop.2.2.1]{Il}), so it suffices
to check the assertion for the ring $B$ and its regular element
$T$. However, the sequence of ring homomorphisms $A\to B\to B_0$
induces a distinguished triangle (\cite[Ch.II, Prop.2.1.2]{Il})
$$
\L_{B/A}\otimes_BB_0\to\L_{B_0/A}\to\L_{B_0/B}\to
\L_{B/A}\otimes_BB_0[1]
$$
and since clearly $\L_{B_0/A}\simeq 0$ in $\sD(B_0\Mod)$,
and $\L_{B/A}\simeq\Omega^1_{B/A}[0]\simeq B[0]$
(\cite[Ch.II, Prop.1.2.4.4]{Il}), the assertion follows
(details left to the reader). Next, suppose that $n>1$,
and that the assertion is already known for regular sequences
of length $<n$. Denote by $I'\subset A$ the ideal generated
by $a_1,\dots,a_{n-1}$, and set $A':=A/I'$. There follows a
sequence of ring homomorphisms $A\to A'\to A_0$, and the
inductive assumption implies that the natural morphisms
$$
\L_{A'/A}\to I'/I'{}^2[1]
\qquad
\L_{A_0/A'}\to I/(I^2+I')[1]
$$
are isomorphisms, and $I'/I'{}^2$ (resp. $I/(I^2+I')$) is a
free $A'$-module (resp. $A_0$-module) of rank $n-1$ (resp.
of rank $1$). Then the assertion follows easily, by
inspecting the distinguished triangle 
$$
\L_{A'/A}\otimes_{A'}A_0\to\L_{A_0/A}\to\L_{A_0/A'}\to
\L_{A'/A}\otimes_{A'}A_0[1]
$$
given again by \cite[Ch.II, Prop.2.1.2]{Il} (details left
to the reader).
\end{proof}

\begin{proposition}\label{prop_cover-n-1}
Let $A$ be a local noetherian ring, $a\in A$ a non-invertible
element, and set $A_0:=A/aA$. The following conditions are
equivalent :
\begin{enumerate}
\alphaenu
\item
$a$ is a regular element of $A$.
\item
$H_2\L_{A_0/A}=0$, and $aA/a^2A$ is a free $A_0$-module of rank one.
\end{enumerate}
\end{proposition}
\begin{proof} For any ring $R$, any non-invertible element
$x\in R$, and any $n\in\N$, set $R_n:=R/x^{n+1}R$, and consider
the $R_0$-linear map
$$
\beta_{x,n}:R_0\to x^nR/x^{n+1}R
$$
induced by multiplication by $x^n$. We remark :

\begin{claim}\label{cl_put-it-first}
Suppose that $R$ is a noetherian local ring, denote by
$\kappa_R$ the residue field of $R$, and let $n>0$
be any given integer. The following conditions
are equivalent :
\begin{enumerate}
\alphaenu
\addenu\addenu
\item
$\beta_{x,n}$ is an isomorphism.
\item
$x^n\neq 0$ and $\Tor_1^{R_0}(x^nR/x^{n+1}R,\kappa_R)=0$.
\item
$x^n\neq 0$ and the surjection $R_0\to\kappa_R$ induces
a surjective map
$$
H_2(\L_{R_{n-1}/R}\otimes_{R_{n-1}}R_0)\to
H_2(\L_{R_{n-1}/R}\otimes_{R_{n-1}}\kappa_R).
$$
\end{enumerate}
\end{claim}
\begin{pfclaim} It is easily seen that (c)$\Rightarrow$(d).

Conversely, if (d) holds, notice that $x^nR/x^{n+1}R\neq 0$,
since $\bigcap_{n\in\N}x^nR=0$.
Hence $\kappa_R\otimes_R\beta_{x,n}$ is an isomorphism of
one-dimensional $\kappa_R$-vector spaces. On the other hand,
under assumption (d), the natural map
$\kappa_R\otimes_R\Ker\,\beta_{x,n}\to
\Ker(\kappa_R\otimes_R\beta_{x,n})$
is an isomorphism. By Nakayama's lemma, we conclude that
$\Ker\,\beta_{x,n}=0$, {\em i.e.} (c) holds.

Next, recall the natural isomorphism of $R_0$-modules
$$
H_1(\L_{R_{n-1}/R}\otimes_{R_{n-1}}M)\isom
x^nR/x^{2n}R\otimes_{R_{n-1}}M\isom x^nR/x^{n+1}R\otimes_{R_0}M
$$
for every $R_0$-module $M$ (\cite[Ch.III, Cor.1.2.8.1]{Il}).
Denote by $\fm_0$ the kernel of the surjection $R_0\to\kappa_R$;
there follows a left exact sequence
$$
0\to\Tor_1^{R_0}(x^nR/x^{n+1}R,\kappa_R)\to
H_1(\L_{R_{n-1}/R}\otimes_{R_{n-1}}\fm_0)\to
H_1(\L_{R_{n-1}/R}\otimes_{R_{n-1}}R_0)
$$
which shows that (d)$\Leftrightarrow$(e) (details left to the reader).
\end{pfclaim}

\begin{claim}\label{cl_reg-crit-Tor}
In the situation of claim \ref{cl_put-it-first}, the following
conditions are equivalent :
\begin{enumerate}
\alphaenu
\addenu\addenu\addenu\addenu\addenu
\item
$x$ is a regular element of $R$.
\item
$\beta_{x,n}$ is an isomorphism, for every $n\in\N$.
\end{enumerate}
\end{claim}
\begin{pfclaim} If (f) holds, $x^n$ is a regular element of $R$
for every $n>0$, and then (g) follows easily. Conversely, assume
(g), and suppose that $yx=0$ for some $y\in R$; we claim that
$y\in x^nR$ for every $n\in\N$. We argue by induction on
$n$ : for $n=0$, there is nothing to prove. Suppose that we
have already obtained a factorization $y=x^nz$ for some
$z\in R$. Then $x^{n+1}z=0$, so the class of $z$ in $R_0$
lies in $\Ker\,\beta_{n+1}$, hence this class must vanish,
{\em i.e.} $z\in xR$, and therefore $y\in x^{n+1}R$. Since
$\bigcap_{n\in\N}x^nR=0$, we deduce that $y=0$, whence (f).
\end{pfclaim}

\begin{claim}\label{cl_go-up}
Let $f:R\to R'$ be any ring homomorphism, and set $x':=f(x)$.
Suppose that $\beta_{x,n}$ and $\beta_{x',n}$ are both
isomorphisms, for some integer $n\in\N$, and set
$R'_n:=R'/x'{}^{n+1}R'$. Then the induced morphism
$$
\L_{R_0/R_n}\otimes_{R_0}R'_0\to\L_{R'_0/R'_n}
$$
is an isomorphism in $\sD(R'_0\Mod)$.
\end{claim}
\begin{pfclaim} If $n=0$, there is nothing to prove, hence
assume that $n>0$. We remark that, for every $i=0,\dots,n$,
the complex
$$
R_n\xrightarrow{\ x^{i+1}\ }R_n\xrightarrow{\ x^{n-i}\ }R_n
$$
is exact. Indeed, for $i=0$, this results immediately from
the assumption that $\Ker\,\beta_{n,x}=0$. Suppose that $i>0$,
and that the assertion is already known for $i-1$; then,
if $yx^{n-i}\in x^{n+1}R$ for some $y\in R$, the inductive
hypothesis yields $y\in x^iR$, so say that $y=x^iu$ and
$yx^{n-i}=zx^{n+1}$ for some $u,z\in R$. It follows that
$x^n(u-zx)=0$, and then $u-zx\in xR$, again since
$\Ker\,\beta_{x,n}=0$; thus, $u\in xR$, and $y\in x^{i+1}R$,
as asserted.

We deduce that the $R_n$-module $R_0$ admits a free resolution
$$
\Sigma
\quad :\quad
\cdots\to R_n\xrightarrow{\ x\ }R_n\xrightarrow{\ x^n\ }R_n
\xrightarrow{\ x\ }R_n\to R_0.
$$
The same argument applies to $R'$ and its element $x'$, and
yields a corresponding free resolution $\Sigma'$ of the
$R'_n$-module $R'_0$. A simple inspection shows that
$\Sigma'\otimes_{R_n}R'_n=\Sigma$, {\em i.e.} the natural
morphism $R_0\derotimes_{R_n}R'_n\to R'_0$ is an isomorphism
in $\sD(R'_0\Mod)$. The claim then follows from
\cite[Ch.II, Prop.2.2.1]{Il}.
\end{pfclaim}

With these preliminaries, we may now return to the situation
of the proposition : first, lemma \ref{lem_regu-seq-cot} says
that (a)$\Rightarrow$(b). For the converse, we shall apply the
criterion of claim \ref{cl_reg-crit-Tor} : namely, we shall
show, by induction on $n$, that $\beta_{a,n}:A_0\to a^nA/a^{n+1}A$
is bijective for every $n\in\N$.

For $n=0$, there is nothing to prove. Assume that $n>0$, and
that the assertion is already known for $n-1$. Let
$\fm\subset A[T]$ be the (unique) maximal ideal containing $T$,
and set $B:=A_\fm$. We let $f:B\to A$ be the map of $A$-algebras
given by the rule : $T\mapsto a$. Define as usual
$B_n:=B/T^{n+1}B$ and $A_n:=A/a^{n+1}A$ for every $n\in\N$.
Clearly $\beta_{T,n-1}:B_0\to T^{n-1}B/T^nB$ is bijective,
and the same holds for $\beta_{a,n-1}$, by inductive assumption.
Then, claim \ref{cl_go-up} says that the induced morphism
$\L_{B_0/B_{n-1}}\otimes_{B_0}A_0\to\L_{A_0/A_{n-1}}$ is an
isomorphism in $\sD(A_0\Mod)$. Denote by $\kappa$ the residue
field of $A$ and $B$, and notice as well that
$H_2(\L_{A_0/A}\otimes_{A_0}\kappa)=0$, by virtue of (b).
Consequently, the commutative diagram of ring homomorphisms
$$
\xymatrix{
B \ar[r] \ar[d] & B_{n-1} \ar[r] \ar[d] & B_0 \ar[d] \\
A \ar[r] & A_{n-1} \ar[r] & A_0
}$$
induces a commutative ladder with exact rows
(\cite[Ch.II, Prop.2.1.2]{Il}) :
$$
\xymatrix{
H_3(\L_{B_0/B_{n-1}}\otimes_{B_0}\kappa) \ar[r] \ar[d] &
H_2(\L_{B_{n-1}/B}\otimes_{B_{n-1}}\kappa) \ar[d] \ar[r] &
H_2(\L_{B_0/B}\otimes_{B_0}\kappa) \ar[d] \\
H_3(\L_{A_0/A_{n-1}}\otimes_{A_0}\kappa) \ar[r] &
H_2(\L_{A_{n-1}/A}\otimes_{A_{n-1}}\kappa) \ar[r] & 0
}$$
whose left vertical arrow is an isomorphism. It follows
that the central vertical arrow is surjective. Consider
now the commutative diagram
\set\begin{equation}\label{eq_switch-to-right}
{\diagram
H_2(\L_{B_{n-1}/B}\otimes_{B_{n-1}}B_0) \ar[r] \ar[d] &
H_2(\L_{B_{n-1}/B}\otimes_{B_{n-1}}\kappa) \ar[d] \\
H_2(\L_{A_{n-1}/A}\otimes_{A_{n-1}}A_0) \ar[r] &
H_2(\L_{A_{n-1}/A}\otimes_{A_{n-1}}\kappa)
\enddiagram}
\end{equation}
induced by the maps $B_0\to A_0\to\kappa$. We have just seen
that the right vertical arrow of \eqref{eq_switch-to-right}
is surjective, and the same holds for its top horizontal
arrow, in light of claim \ref{cl_put-it-first}. Thus, finally,
the bottom horizontal arrow is surjective as well, so
$\beta_{a,n}$ is an isomorphism (claim \ref{cl_put-it-first}),
and the proposition is proved.
\end{proof}

\begin{theorem}\label{th_Andre-reg-crit}
Let $A$ be a noetherian local ring, $I\subset A$ an ideal,
and set $A_0:=A/I$. The following conditions are equivalent :
\begin{enumerate}
\alphaenu
\item
Every minimal system of generators of $I$ is a regular
sequence of elements of $A$.
\item
$I$ is generated by a regular sequence of $A$.
\item
The natural morphism $\L_{A_0/A}\to I/I^2[1]$ is an isomorphism
in $\sD(A_0\Mod)$, and $I/I^2$ is a flat $A_0$-module.
\item
$H_2\L_{A_0/A}=0$, and $I/I^2$ is a flat $A_0$-module.
\end{enumerate}
\end{theorem}
\begin{proof} Clearly (a)$\Rightarrow$(b) and (c)$\Rightarrow$(d);
also, lemma \ref{lem_regu-seq-cot} shows that (b)$\Rightarrow$(c).

(d)$\Rightarrow$(a): Let $(a_1,\dots,a_n)$ be a minimal
system of generators for $I$; recall that the length $n$
of the sequence equals $\dim_\kappa I\otimes_A\kappa$, where
$\kappa$ denotes the residue field of $A$. We shall argue by
induction on $n$. For $n=0$, there is nothing to show, and
the case $n=1$ is covered by proposition \ref{prop_cover-n-1}.
Set $B:=A/a_1A$ and $J:=IB$. Assumption (d) implies that
$H_2(\L_{A_0/A}\otimes_{A_0}\kappa)=0$, therefore the sequence
of ring homomorphisms $A\to B\to A_0$ induces an exact sequence
$$
0\to H_2(\L_{A_0/B}\otimes_{A_0}\kappa)\to
H_1(\L_{B/A}\otimes_B\kappa)\to I\otimes_A\kappa
\to J\otimes_B\kappa\to 0
$$
(\cite[Ch.II, Prop.2.1.2 and Ch.III, Cor.1.2.8.1]{Il}).
However, clearly $J$ admits a generating system of length
$n-1$, hence $n':=\dim_\kappa J\otimes_B\kappa<n$. On the
other hand, $\dim_\kappa H_1(\L_{B/A}\otimes_B\kappa)=1$, so
we have necessarily $n'=n-1$ and
$H_2(\L_{A_0/B}\otimes_{A_0}\kappa)=0$. The latter means
that $H_2\L_{A_0/B}=0$ and $J/J^2$ is a flat $B$-module.
By inductive assumption, we deduce that the sequence
$(\bar a_2,\dots,\bar a_n)$ of the images in $B$ of
$(a_2,\dots,a_n)$, is regular. By virtue of lemma
\ref{lem_regu-seq-cot}, it follows that
$H_3(\L_{A_0/B}\otimes_{A_0}\kappa)=0$, whence a left exact
sequence
$$
0\to H_2(\L_{B/A}\otimes_B\kappa)\to
H_2(\L_{A_0/A}\otimes_{A_0}\kappa)=0
$$
obtained by applying again
\cite[Ch.II, Prop.2.1.2 and Ch.III, Cor.1.2.81]{Il} to the
sequence $A\to B\to A_0$. Thus,
$H_2\L_{B/A}=0$ and $a_1A/a_1^2A$ is a flat $B$-module, so
$a_1$ is a regular element (proposition \ref{prop_cover-n-1})
and finally, $(a_1,\dots,a_n)$ is a regular sequence, as
required.
\end{proof}

\begin{corollary}\label{cor_Andre-reg-crit}
Let $A$ be a local noetherian ring, with maximal
ideal $\fm$, and residue field $\kappa$. Then the following
conditions are equivalent :
\begin{enumerate}
\alphaenu
\item
$A$ is regular.
\item
The natural morphism $\L_{\kappa/A}\to\fm/\fm^2[1]$ is an
isomorphism.
\item
$H_2\L_{\kappa/A}=0$.
\end{enumerate}
\end{corollary}
\begin{proof} It follows immediately, by invoking theorem
\ref{th_Andre-reg-crit} with $I:=\fm$, and lemma
\ref{lem_regu-seq-cot}.
\end{proof}

\subsection{Excellent rings}\label{sec_excellent-rings}
Recall that a morphism of schemes $f:X\to Y$ is called
{\em regular}, if it is flat, and for every $y\in Y$, the
fibre $f^{-1}(y)$ is locally noetherian and regular
(\cite[Ch.IV, D{\'e}f.6.8.1]{EGAIV-2}).

\begin{lemma}\label{lem_yoga-reg}
Let $f:X\to Y$ and $g:Y\to Z$ be two morphisms of locally
noetherian schemes.
We have :
\begin{enumerate}
\item
If $f$ and $g$ are regular, then the same holds for $g\circ f$.
\item
If $g\circ f$ is regular, and $f$ is faithfully flat, then $g$
is regular.
\end{enumerate}
\end{lemma}
\begin{proof}(i): Clearly $h:=g\circ f$ is flat. Let $z\in Z$
be any point, and $K$ any finite extension of $\kappa(z)$. Set
$$
X':=h^{-1}(z)\times_{\kappa(z)}K
\qquad\text{and}\qquad
Y':=g^{-1}(z)\times_{\kappa(z)}K.
$$
It is easily seen that the induced morphism $f':X'\to Y'$
is regular. Moreover, for every $y\in Y'$, the local ring
$\cO_{Y',y'}$ is regular, since $g$ is regular. Then the
assertion follows from :

\begin{claim} Let $A\to B$ be a flat and local ring homomorphism
of local noetherian rings. Denote by $\fm_A\subset A$ the
maximal ideal, and suppose that both $A$ and $B_0:=B/\fm_AB$
are regular. Then $B$ is regular.
\end{claim}
\begin{pfclaim} On the one hand, $\dim B=\dim A+\dim B_0$
(\cite[Ch.IV, Cor.6.1.2]{EGAIV}).
On the other hand, let $a_1,\dots,a_n$ be a minimal generating
system for $\fm_A$, and $b_1,\dots,b_m$ a system of elements
of the maximal ideal $\fm_B$ of $B$, whose images in $B_0$ is
a minimal generating system for $\fm_B/\fm_AB$. By Nakayama's
lemma, it is easily seen that the system
$a_1,\dots,a_n,b_1,\dots,b_m$ generates the ideal $\fm_B$.
Since $A$ and $B_0$ are regular, $n=\dim A$ and $m=\dim B$,
so $n+m=\dim B$, and the claim follows.
\end{pfclaim}

(ii): Clearly $g$ is flat. Then the assertion follows easily
from \cite[Ch.0, Prop.17.3.3(i)]{EGAIV} : details left to the
reader.
\end{proof}

\begin{definition}\label{def_excelll}
Let $A$ be a noetherian ring.
\begin{enumerate}
\item
We say that $A$ is a {\em G-ring}, if the formal fibres of
$\Spec\,A$ are geometrically regular, {\em i.e.} for every
$\fp\in\Spec\,A$, the natural morphism
$\Spec\,A^\wedge_\fp\to\Spec\,A_\fp$ from the spectrum of the
$\fp$-adic completion of $A$, is regular : see
\cite[Ch.IV, \S7.3.13]{EGAIV-2}.
\item
We say that $A$ is {\em quasi-excellent}, if $A$ is a
G-ring, and moreover the following holds.
For every prime ideal $\fp\subset A$, and every finite
radicial extension $K'$ of the field of fractions $K$ of
$B:=A/\fp$, there exists a finite $B$-subalgebra $B'$
of $K'$ such that the field of fractions of $B'$ is $K'$,
and the {\em regular locus\/} of $\Spec\,B'$ is an open subset
(the latter is the set of all prime ideals $\fq\subset B'$
such that $B'_\fq$ is a regular ring).
\item
We say that $A$ is a {\em Nagata ring}, if the following
holds. For every $\fp\in\Spec\,A$ and every finite
field extension $\kappa(\fp)\subset L$, the integral
closure of $A/\fp$ in $L$ is a finite $A$-module. (The rings
enjoying this latter property are called {\em universally
japanese\/} in \cite[Ch.0, D\'ef.23.1.1]{EGAIV}.)
\item
We say that $A$ is {\em universally catenarian\/} if every
$A$-algebra $B$ of finite type is catenarian, {\em i.e.} any
two saturated chains $(\fp_0\subset\cdots\subset\fp_n)$,
$(\fq_0\subset\cdots\subset\fq_m)$ of prime ideals of $B$,
with $\fp_0=\fq_0$ and $\fp_n=\fq_m$, have the same length
(so $n=m$) (\cite[Ch.IV, D{\'e}f.5.6.2]{EGAIV-2}).
\item
We say that $A$ is {\em excellent}, if it is quasi-excellent
and universally catenarian (\cite[Ch.IV, D{\'e}f.7.8.2]{EGAIV-2}).
\end{enumerate}
\end{definition}

\begin{lemma}\label{lem_Nagatas}
Let $A$ be a noetherian ring.
\begin{enumerate}
\item
If $A$ is quasi-excellent, then $A$ is a Nagata ring.
\item
If $A$ is a G-ring, then every quotient and every localization
of $A$ is a G-ring.
\item
Suppose that the natural morphism $\Spec\,A^\wedge_\fm\to\Spec\,A_\fm$
is regular for every maximal ideal $\fm\subset A$. Then $A$ is
a G-ring.
\item
If $A$ is a local G-ring, then $A$ is quasi-excellent.
\item
If $A$ is a complete local ring, then $A$ is excellent.
\end{enumerate}
\end{lemma}
\begin{proof}(i): This is \cite[Ch.IV, Cor.7.7.3]{EGAIV-2}.

(ii): The assertion for localizations is obvious. Next, if
$I\subset A$ is any ideal, and $\fp\subset A$ any prime ideal
containing $I$, then
$(A/I)^\wedge_\fp=A^\wedge_\fp/IA^\wedge_\fp$, from which it
is immediate that $A/I$ is a G-ring, if the same holds for $A$.

(iv) follows from
\cite[Ch.IV, Th.6.12.7, Prop.7.3.18, Th.7.4.4(ii)]{EGAIV-2}.

(v): In light of (iv), it suffices to remark that every complete
noetherian local ring is universally catenarian
(\cite[Ch.IV, Prop.5.6.4]{EGAIV-2} and
\cite[Ch.0, Th.19.8.8(i)]{EGAIV}) and is a G-ring
(\cite[Ch.0, Th.22.3.3, Th.22.5.8, and Prop.19.3.5(iii)]{EGAIV}).

(iii): Let us remark, more generally :

\begin{claim}\label{cl_descent-G-ring}
Let $\phi:A\to B$ be a faithfully flat ring homomorphism of
noetherian rings, such that $f:=\Spec\,\phi$ is regular. If
$B$ is a G-ring, the same holds for $A$.
\end{claim}
\begin{pfclaim} In light of (ii), we easily reduce to the case
where both $A$ and $B$ are local, $\phi$ is a local ring
homomorphism, and it suffices to show that the natural
morphism $\pi_A:\Spec\,A^\wedge\to\Spec\,A$ is regular
(where $A^\wedge$ is the completion of $A$).
Consider the commutative diagram :
\set\begin{equation}\label{eq_f-and-f-wedge}
{\diagram \Spec\,B^\wedge \ar[r]^-{f^\wedge} \ar[d]_{\pi_B} &
          \Spec\,A^\wedge \ar[d]^{\pi_A} \\
          \Spec\,B \ar[r]^-f & \Spec\,A.
\enddiagram}
\end{equation}
By assumption, $\pi_B$ is a regular morphism; Then the same
holds for $f\circ\pi_B=\pi_A\circ f^\wedge$
(lemma \ref{lem_yoga-reg}(i)). However, it is easily seen
that the induced map $\phi^\wedge:A^\wedge\to B^\wedge$
is still a local ring homomorphism, hence $f^\wedge$ is
faithfully flat, so the claim follows from lemma
\ref{lem_yoga-reg}(ii).
\end{pfclaim}

Now, in order to prove (iii), it suffices to check that
$A_\fm$ is a G-ring for every maximal ideal $\fm\subset A$.
In view of our assumption, the latter assertion follows
from claim \ref{cl_descent-G-ring} and (v).
\end{proof}

\begin{proposition}\label{prop_vanish-cot-n}
Let $A$ be a ring, $B$ a noetherian $A$-algebra of finite
(Krull) dimension, $n\in\N$ an integer, and suppose that
$H_k(\L_{B/A}\otimes_B\kappa(\fp))=0$ for every prime ideal
$\fp\subset B$ and every $k=n,\dots,n+\dim B$.
Then $H_n(\L_{B/A}\otimes_BM)=0$ for every $B$-module $M$.
\end{proposition}
\begin{proof} Let us start out with the following more general :

\begin{claim}\label{cl_primes-cot}
Let $A$ be a ring, $B$ a noetherian $A$-algebra, $\fp\in\Spec\,B$
a prime ideal, $n\in\N$ an integer, and suppose that the following
two conditions hold :
\begin{enumerate}
\alphaenu
\item
$H_n(\L_{B/A}\otimes_B\kappa(\fp))=0$.
\item
$H_{n+1}(\L_{B/A}\otimes_BB/\fq)=0$ for every proper specialization
$\fq$ of $\fp$ in $\Spec\,B$.
\end{enumerate}
Then, $H_n(\L_{B/A}\otimes_BB/\fp)=0$.
\end{claim}
\begin{pfclaim} Let $b\in B\setminus\fp$ be any element, and set
$M:=B/(\fp+bB)$. Since $B$ is noetherian, $M$ admits a finite
filtration $M_0\subset M_1\subset\cdots\subset M_n:=M$ such
that, for every $i=0,\dots,n-1$, the subquotient $M_{i+1}/M_i$
is isomorphic to $B/\fq$, for some proper specialization $\fq$
of $\fp$ (\cite[Th.6.4]{Mat}). From (b), and a simple induction,
we deduce that $H_{n+1}(\L_{B/A}\otimes_BM)=0$. Whence, by
considering the short exact sequence of $B$-modules
$$
0\to B/\fp\xrightarrow{\ b\ }B/\fp\to M\to 0
$$
we see that scalar multiplication by $b$ is an injective map
on the $B$-module $H_n(\L_{B/A}\otimes_BB/\fp)$. Since this
holds for every $b\in B\setminus\fp$, we conclude that the
natural map
$$
H_n(\L_{B/A}\otimes_BB/\fp)\to
H_n(\L_{B/A}\otimes_BB/\fp)\otimes_BB_\fp=
H_n(\L_{B/A}\otimes_B\kappa(\fp))
$$
is injective. Then the assertion follows from (a).
\end{pfclaim}

Now, let $\fp\subset B$ be any prime ideal; the assumption,
together with claim \ref{cl_primes-cot} and a simple induction
on $d:=\dim B/\fp$, shows that
$$
H_k(\L_{B/A}\otimes_BB/\fp)=0
\qquad
\text{for every $k=n,\dots,n+\dim B-d$.}
$$
For any $B$-module $M$, set $H(M):=H_n(\L_{B/A}\otimes_BM)$.
Especially, we get $H(B/\fp)=0$, for any prime ideal
$\fp\subset B$. Since $B$ is noetherian, it follows easily that
$H(M)=0$ for any $B$-module $M$ of finite type (details left to
the reader). Next, if $M$ is arbitrary, we may write it as the
union of the filtered family $(M_i~|~i\in I)$ of its submodules
of finite type; since $H(M)$ is the colimit of the induced system
$(H(M_i)~|~i\in I)$, we see that $H(M)=0$, as sought.
\end{proof}

\begin{corollary}\label{cor_reg-maps-crit}
Let $A\to B$ be a homomorphism of noetherian rings. Then
the following conditions are equivalent :
\begin{enumerate}
\alphaenu
\item
$\Omega^1_{B/A}$ is a flat $B$-module, and $H_i\L_{B/A}=0$ for every
$i>0$.
\item
$\Omega^1_{B/A}$ is a flat $B$-module, and $H_1\L_{B/A}=0$.
\item
The induced morphism of schemes $f:\Spec\,B\to\Spec\,A$ is regular.
\item
$H_1(\L_{B/A}\otimes_B\kappa(x))=0$ for every $x\in\Spec\,B$.
\end{enumerate}
\end{corollary}
\begin{proof} Let $x\in X:=\Spec\,B$ be any point; according to
\cite[Ch.0, Th.19.7.1]{EGAIV}, the following conditions are
equivalent :
\begin{enumerate}
\alphaenu
\addenu\addenu\addenu\addenu
\item
the map on stalks $\cO_{Y,f(x)}\to\cO_{\!X,x}$ is formally
smooth for the preadic topologies defined by the maximal ideals.
\item
$f$ is flat at the point $x$, and the $\kappa(f(x))$-algebra
$\cO_{f^{-1}f(x),x}$ is geometrically regular.
\end{enumerate}
On the other hand, by virtue of proposition \ref{prop_fsmooth-cot},
condition (e) is equivalent to the vanishing of
$H_1(\L_{B/A}\otimes_B\kappa(x))$, whence
(b)$\Rightarrow$(c)$\Leftrightarrow$(d). It remains to check that
(d)$\Rightarrow$(a). However, assume (d); taking into account
proposition \ref{prop_vanish-cot-n}, and arguing by induction
on $i$, we easily show that (a) will follow, provided
$$
H_i(\L_{B/A}\otimes_B\kappa(x))=0
\qquad
\text{for every $x\in X$ and every $i>1$}.
$$
For every $x\in X$, set $B_x:=\cO_{\!X,x}\otimes_A\kappa(f(x))$;
since $B$ is $A$-flat, the latter holds if and only if
$$
H_i(\L_{B_x/\kappa(f(x))}\otimes_B\kappa(x))=0
\qquad
\text{for every $x\in X$ and every $i>1$}
$$
(\cite[Ch.II, Prop.2.2.1 and Ch.III, Cor.2.3.1.1]{Il}).
Hence, we may replace $A$ by $\kappa(f(x))$, and $B$ by $B_x$,
and assume from start that $A$ is a field and $B$ is a local
geometrically regular $A$-algebra with residue field $\kappa_B$,
and it remains to check that $H_i(\L_{B/A}\otimes_B\kappa_B)=0$
for every $i>1$. However, in view of corollary \ref{cor_Andre-reg-crit},
the sequence of ring homomorphisms $A\to B\to\kappa_B$ yields
an isomorphism
$$
H_i(\L_{B/A}\otimes_B\kappa_B)\isom H_i\L_{\kappa_B/A}
\qquad
\text{for every $i>1$}
$$
(\cite[Ch.II, Prop.2.1.2]{Il}) so we conclude by the following general :

\begin{claim} Let $K\subset E$ be any extension of fields.
Then $H_i\L_{E/K}=0$ for every $i>1$.
\end{claim}
\begin{pfclaim}[] By \cite[Ch.II, (1.2.3.4)]{Il}, we may reduce
to the case where $E$ is a finitely generated extension of $K$,
say $E=K(a_1,\dots,a_n)$. We proceed by induction on $n$. If
$n=1$, then $E$ is either an algebraic extension or a purely
transcendental extension of $K$. In the latter case, the
assertion is immediate
(\cite[Ch.II, Prop.1.2.4.4 and Ch.III, Cor.2.3.1.1]{Il}).
For the case of an algebraic extension, the assertion is a special
case of \cite[Th.6.3.32(i)]{Ga-Ra} (we apply {\em loc.cit.} to the
valued field $(K,|\cdot|)$ with trivial valuation $|\cdot|$). Lastly,
if $n>1$, set $L:=K(a_1,\dots,a_{n-1})$. By inductive assumption we
have $H_i\L_{L/K}=H_i\L_{E/L}=0$ for $i>1$; on the other hand,
there is a distinguished triangle (\cite[Ch.II, Prop.2.1.2]{Il})
$$
\L_{L/K}\otimes_LE\to\L_{E/K}\to\L_{E/L}\to\L_{L/K}\otimes_LE[1]
\qquad
\text{in $\sD(E\Mod)$}.
$$
The sought vanishing follows immediately, 
\end{pfclaim}
\end{proof}

\begin{corollary}\label{cor_crit-q.excel}
Let $A$ be a local noetherian ring, $A^\wedge$ the completion
of $A$. Then the following conditions are equivalent :
\begin{enumerate}
\alphaenu
\item
$A$ is quasi-excellent.
\item
$\Omega^1_{A^\wedge/A}$ is a flat $A^\wedge$-module, and
$H_1\L_{A^\wedge/A}=0$.
\item
$\Omega^1_{A^\wedge/A}$ is a flat $A^\wedge$-module, and
$H_i\L_{A^\wedge/A}=0$ for every $i>0$.
\end{enumerate}
\end{corollary}
\begin{proof} Taking into account lemma \ref{lem_Nagatas}(iii),
this is a special case of corollary \ref{cor_reg-maps-crit}.
\end{proof}

\begin{proposition}\label{prop_separation-degree}
Let $p>0$ be a prime integer, $A$ a regular local and excellent
$\F_p$-algebra, $B$ a local noetherian $A$-algebra, $\fm_B$ the
maximal ideal of $B$, and $M$ a $B$-module of finite type. Then
the $B$-module $\Omega^1_{A/\F_p}\otimes_AM$ is separated for
the $\fm_B$-preadic topology.
\end{proposition}
\begin{proof} Let $\phi:A\to B$ be the structure morphism,
and set $\fp:=\phi^{-1}\fm_B$; the localization $A_\fp$ is
still excellent (lemma \ref{lem_Nagatas}(ii,iv)) and regular,
and clearly $\Omega^1_{A_\fp/\F_p}\otimes_{A_\fp}M=
\Omega^1_{A/\F_p}\otimes_AM$, hence we may replace $A$ by
$A_\fp$, and assume that $\phi$ is local. Let $A^\wedge$
and $B^\wedge$ be the completions of $A$ and $B$, set
$M^\wedge:=B^\wedge\otimes_BM$, and notice that both of the
natural maps $\F_p\to A$ and $A\to A^\wedge$ are regular.
In view of corollary \ref{cor_reg-maps-crit}, it follows
that both of the natural $B$-linear maps
$$
\Omega^1_{A/\F_p}\otimes_AM\to\Omega^1_{A/\F_p}\otimes_AM^\wedge
\qquad
\Omega^1_{A/\F_p}\otimes_AM\to
\Omega^1_{A^\wedge/\F_p}\otimes_{A^\wedge}M^\wedge
$$
are injective (to see the injectivity of the second map,
one applies the transitivity triangle arising from the
sequence of ring homomorphisms $\F_p\to A\to A^\wedge$ :
details left to the reader). Thus, we may assume that
$A$ is complete. Now, for every ring $R$, and every
integer $m\in\N$, set
\set\begin{equation}\label{eq_give-me-flat}
R_{(m)}:=R[[T_1,\dots,T_m]]
\qquad
R_{\La m\Ra}:=R[[T_1^p,\dots,T^p_m]]\subset R_{(m)}.
\end{equation}
With this notation, we have an isomorphism
$A\isom\kappa_{(d)}$ of $\F_p$-algebras, where $\kappa$
is the residue field of $A$, and $d:=\dim A$
(\cite[Ch.0, Th.19.6.4]{EGAIV}).

\begin{claim}\label{cl_general-field}
Let $K$ be any field of characteristic $p$, and $m\in\N$ any integer.
We have :
\begin{enumerate}
\item
There exists a cofiltered system $(K^\lambda~|~\lambda\in\Lambda)$
of subfields of $K$ such that $[K:K^\lambda]$ is finite for every
$\lambda\in\Lambda$, and $\bigcap_{\lambda\in\Lambda}K^\lambda=K^p$.
\item
For every system $(K^\lambda~|~\lambda\in\Lambda)$ fulfilling the
condition of (i), the following holds :
\begin{enumerate}
\item
$\Omega^1_{K_{(m)}/K^\lambda_{\La m\Ra}}$ is a free $K_{(m)}$-module
of finite rank, for every $\lambda\in\Lambda$, and the rule
$$
M\mapsto\Omega_m(M):=\lim_{\lambda\in\Lambda}\,
(\Omega^1_{K_{(m)}/K^\lambda_{\La m\Ra}}\otimes_{K_{(m)}}M)
$$
defines an exact functor $K_{(m)}\Mod\to K_{(m)}\Mod$.
\item
The natural map
$$
\eta_m(M):\Omega^1_{K_{(m)}/\F_p}\otimes_{K_{(m)}}M\to\Omega_m(M)
$$
is injective for every $K_{(m)}$-module $M$.
\item
Let $F$ (resp. $F^\lambda$) denote the field of fractions of
$K_{(m)}$ (resp. of $K^\lambda_{\La m\Ra}$, for every
$\lambda\in\Lambda$); then $\bigcap_{\lambda\in\Lambda}F^\lambda=F^p$.
\end{enumerate}
\end{enumerate}
\end{claim}
\begin{pfclaim} (i) and (ii.c) follow from
\cite[Ch.0, Prop.21.8.8]{EGAIV} (and its proof).

(ii.a): For given $\lambda\in\Lambda$, say that $x_1,\dots,x_r$
is a $p$-basis of $K$ over $K^\lambda$; then it is easily seen
that $x_1,\dots,x_r,T_1,\dots,T_m$ is a $p$-basis of $K_{(m)}$
over $K^\lambda_{\La m\Ra}$ (see \cite[Ch.0, D\'ef.21.1.9]{EGAIV}).
According to \cite[Ch.0, Cor.21.2.5]{EGAIV}, it follows that
$\Omega^1_{K_{(m)}/K^\lambda_{\La m\Ra}}$ is the free
$K_{(m)}$-module of finite type with basis
$dx_1,\dots,dx_r,dT_1,\dots,dT_m$. Moreover, say that
$K^\mu\subset K^\lambda$, and let $x_{r+1},\dots,x_s$ be
a $p$-basis of $K^\lambda$ over $K^\mu$; then $x_1,\dots,x_s$
is a $p$-basis of $K$ over $K^\mu$ (\cite[Ch.0, Lemme 21.1.10]{EGAIV}),
so the induced map
$$
\Omega^1_{K_{(m)}/K^\mu_{\La m\Ra}}\to
\Omega^1_{K_{(m)}/K^\lambda_{\La m\Ra}}
$$
is a projection onto a direct factor, and the assertion follows
easily.

(ii.b): We are easily reduced to the case where $M$ is a
$K_{(m)}$-module of finite type, and in light of (ii.a), we
may further assume that $M$ is a cyclic $K_{(m)}$-module. Next,
we remark that, due to (ii.c), the natural map
$$
\Omega^1_{F/\F_p}\to\lim_{\lambda\in\Lambda}\Omega^1_{F/F^\lambda}
$$
is injective (\cite[Ch.0, Th.21.8.3]{EGAIV}); in other words,
$\eta_m(F)$ is injective. In order to show the injectivity of
$\eta_m(M)$, it then suffices to check that the functor
$M\mapsto\Omega^1_{K_{(m)}/\F_p}\otimes_{K_{(m)}}M$ is exact.
The latter holds by virtue of corollary \ref{cor_reg-maps-crit},
since the (unique) morphism of schemes $\Spec\,K_{(m)}\to\Spec\,\F_p$
is obviously regular. This completes the proof for $m=0$.
Suppose now that $m>0$, and that the injectivity of $\eta_n(M)$
is already known for every $n<m$ and every $K_{(n)}$-module
$M$. By the foregoing, it remains to check that
$\eta_m(K_{(m)}/I)$ is injective, for every non-zero ideal
$I\subset K_{(m)}$. Pick any non-zero $f\in I$, and set
$R:=K_{(m-1)}$; according to \cite[Ch.VII, n.7, Lemme 3]{BouAC}
and \cite[Ch.VII, n.8, Prop.6]{BouAC}, there exist an automorphism
$\sigma$ of the ring $K_{(m)}$, and elements $g\in R[T_m]$,
$u\in K_{(m)}^\times$ such that $\sigma(f)=u\cdot g$, and
$g=T_m^d+a_1T_m^{d-1}+\cdots+a_d$ for some $d\geq 0$ and certain
elements $a_1,\dots,a_d$ of the maximal ideal of $R$.
However, set $M':=K_{(m)}/\sigma(I)$; in view of the commutative
diagram of $\F_p$-modules :
$$
\xymatrix{ \Omega^1_{K_{(m)}/\F_p}\otimes_{K_{(m)}}M
\ar[rr]^-{\eta_m(M)} \ar[d] & & \Omega_m(M) \ar[d] \\
\Omega^1_{K_{(m)}/\F_p}\otimes_{K_{(m)}}M'
\ar[rr]^-{\eta_m(M')} & & \Omega_m(M')
}$$
(whose vertical arrows are induced by $\sigma$) we see that
$\eta_m(M)$ is injective, if and only if the same holds for
$\eta_m(M')$. Hence, we may replace $I$ by $\sigma(I)$, and
assume that $g\in I$. In this case, set also
$R^\lambda:=K^\lambda_{\La m-1\Ra}$ for every $\lambda\in\Lambda$,
and notice that the natural maps
$$
R[T_m]/g^pR[T_m]\to K_{(m)}/g^pK_{(m)}
\qquad
R^\lambda[T^p_m]/g^pR^\lambda[T^p_m]\to
K^\lambda_{\La m\Ra}/g^pK^\lambda_{\La m\Ra}
$$
are bijective; there follows a commutative diagram of
$K_{(m)}$-modules :
$$
\xymatrix{
\Omega^1_{R[T_m]/\F_p}\otimes_{R[T_m]}M \ar[r]
\ar[d]_{\alpha^\lambda\otimes_{R[T_m]} M} &
\Omega^1_{K_{(m)}/\F_p}\otimes_{K_{(m)}}M
\ar[d]^{\eta_m^\lambda\otimes_{K_{(m)}}M} \\
\Omega^1_{R[T_m]/R^\lambda}\otimes_{R[T_m]}M \ar[r] &
\Omega^1_{K_{(m)}/K^\lambda_{\La m\Ra}}\otimes_{K_{(m)}}M
}$$
whose horizontal arrows are isomorphisms, and
$\eta_m(M)=\lim_{\lambda\in\Lambda}\eta_m^\lambda\otimes_{K_{(m)}}M$.
On the other hand, for every $\lambda\in\Lambda$ we have a
commutative ladder of $R[T_m]$-modules with exact rows :
$$
\Sigma_\lambda
\quad :\quad
{\diagram 0 \ar[r] &
\Omega^1_{R/\F_p}\otimes_RR[T_m] \ar[r]
\ar[d]_{\eta^\lambda_{m-1}\otimes_RR[T_m]} &
\Omega^1_{R[T_m]/\F_p} \ar[r] \ar[d]_{\alpha^\lambda} &
\Omega^1_{R[T_m]/R} \ar[r] \ddouble & 0 \\
0 \ar[r] & \Omega^1_{R/R^\lambda}\otimes_RR[T_m] \ar[r] &
\Omega^1_{R[T_m]/R^\lambda} \ar[r] & \Omega^1_{R[T_m]/R} \ar[r] & 0 
\enddiagram}
$$
and notice that the rows of $\Sigma_\lambda\otimes_{R[T_m]}M$
are still short exact, for every $\lambda\in\Lambda$. By
inductive assumption,
$\lim_{\lambda\in\Lambda}\eta^\lambda_{m-1}\otimes_RM$ is an
injective map; we deduce that the same holds for
$\lim_{\lambda\in\Lambda}\alpha^\lambda\otimes_{R[T]}M$, and
the claim follows.
\end{pfclaim}

Take $K:=\kappa$, and pick any cofiltered system
$(K^\lambda~|~\lambda\in\Lambda)$ as provided by claim
\ref{cl_general-field}(i); in light of claim
\ref{cl_general-field}(ii.b), it now suffices to show that
the resulting $\Omega_d(M)$ is a separated $B$-module, for
every noetherian $\kappa_{(d)}$-algebra $B$ and every $B$-module
$M$ of finite type. However, $\Omega_d(M)$ is a submodule of
$\prod_{\lambda\in\Lambda}(\Omega^1_{K_{(m)}/K^\lambda_{\La m\Ra}}
\otimes_{K_{(m)}}M)$, hence we are reduced to checking that each
direct factor of the latter $B$-module is separated. But
in view of claim \ref{cl_general-field}(ii.a), we see that
each such factor is a finite direct sum of copies of $M$, so
finally we come down to the assertion that $M$ is separated
for the $\fm_B$-adic topology, which is well known.
\end{proof}

\begin{theorem}\label{th_block-buster}
Let $\phi:A\to B$ be a local ring homomorphism of local noetherian
rings. Suppose that $A$ is quasi-excellent, and $\phi$ is formally
smooth for the preadic topologies defined by the maximal ideals.
Then $\Spec\,\phi$ is regular.
\end{theorem}
\begin{proof} This is the main result of \cite{AnII}. We begin
with the following general remark :

\begin{claim}\label{cl_general-remark}
Let $R$ be a local noetherian ring, $\fm_R\subset R$ the maximal
ideal, $H:R\Mod\to R\Mod$ an additive functor, and suppose that
\begin{enumerate}
\alphaenu
\item
$H$ is $R$-linear, {\em i.e.} $H(t\cdot\one_M)=t\cdot H(\one_M)$
for every $R$-module $M$, and every $t\in R$.
\item
$H$ is {\em semi-exact}, {\em i.e.} for any short exact
sequence $0\to M'\to M\to M''\to 0$ of $R$-modules, the
induced sequence $H(M')\to H(M)\to H(M'')$ is exact.
\item
$H$ commutes with filtered colimits.
\item
$H(M)$ is separated for the $\fm_R$-adic topology, for
every $R$-module $M$ of finite type.
\item
$H(R/\fm_R)=0$.
\end{enumerate}
Then $H(M)=0$ for every $R$-module $M$.
\end{claim}
\begin{pfclaim} Since $H$ commutes with filtered colimits,
it suffices to show that $H(M)=0$ for every finitely generated
$R$-module $M$, and since $H$ is semi-exact, a simple induction
reduces further to the case where $M$ is a cyclic $R$-module.
Let now $\cF$ be the family of all ideals $I$ of $R$ such that
$H(R/I)\neq 0$, and suppose, by way of contradiction, that
$\cF\neq\emptyset$; pick a maximal element $J$ of $\cF$, and
set $M:=R/J$. By assumption, $J\neq\fm_R$; thus, let $t\in R$
be a non-invertible element with $t\notin J$; we get an exact
sequence
$$
H(M)\xrightarrow{\ t\ }H(M)\to H(B_{(m+n)}/(J+tB_{(m+n)}))
$$
whose third term vanishes, by the maximality of $J$. On the other
hand, $\bigcap_{n\in\N}t^nH(M)=0$, since $H(M)$ is separated. Thus
$H(M)=0$, contradicting the choice of $J$, and the claim follows.
\end{pfclaim}

Denote by $\kappa$ the residue field of $A$, and for every $m,n\in\N$,
let $\phi_{m,n}:A_{(m)}\to B_{(m+n)}$ be the composition of
$\phi_{(m)}:A_{(m)}\to B_{(m)}$ (the $T$-adic completion of
$\phi\otimes_AA[T_1,\dots,T_m]$) with the natural inclusion map
$B_{(m)}\to B_{(m+n)}$ (notation of \eqref{eq_give-me-flat}).

\begin{claim}\label{cl_m-and-n-reg}
In the situation of the theorem, suppose furthermore that $A$
is either a field or a complete discrete valuation ring of
mixed characteristic. Then, the morphism $\Spec\,\phi_{m,n}$
is regular for every $m,n\in\N$.
\end{claim}
\begin{pfclaim} Set
$$
H(M):=H_1(\L_{B_{(m+n)}/A_{(m)}}\otimes_{B_{(m+n)}}M).
$$
We shall consider separately three different cases :

$\bullet$\ \
Suppose first that $A$ is a field of characteristic $p>0$.
According to corollary \ref{cor_reg-maps-crit}, it suffices
to show that $H(M)$ vanishes for every $B_{(m+n)}$-module $M$.
Notice that the natural map $\F_p\to B_{(m+n)}$ is regular; from the
distinguished triangle (\cite[Ch.II, Prop.2.1.2]{Il})
$$
\L_{A_{(m)}/\F_p}\otimes_{A_{(m)}}B_{(m+n)}\to\L_{B_{(m+n)}/\F_p}\to
\L_{B_{(m+n)}/A_{(m)}}\to\L_{A_{(m)}/\F_p}\otimes_{A_{(m)}}B_{(m+n)}[1]
$$
and corollary \ref{cor_reg-maps-crit}, we deduce an injective
$B_{(m+n)}$-linear map
$$
H(M)\to\Omega^1_{A_{(m)}/\F_p}\otimes_{A_{(m)}}M
$$
from which it follows that if $M$ is a $B_{(m+n)}$-module of finite
type, $H(M)$ is a separated $B_{(m+n)}$-module, for the preadic
topology defined by the maximal ideal of $B_{(m+n)}$ (proposition
\ref{prop_separation-degree}). Since $\phi$
is formally smooth, $\phi_{m,n}$ is also formally smooth for the
preadic topologies defined by the maximal ideals of $A_{(m)}$ and
$B_{(m+n)}$; from proposition \ref{prop_fsmooth-cot}, we see that
$H(M)=0$, if $M$ is the residue field of $B_{(m+n)}$. Then the
assertion follows from claim \ref{cl_general-remark}.

$\bullet$\ \ 
Next, suppose that $A$ is either a field of characteristic zero,
or a complete discrete valuation ring of mixed characteristic
(so either $\kappa=A$, or else $\kappa$ is a field of positive
characteristic). According to corollary \ref{cor_reg-maps-crit},
it suffices to show that $H(M)$ vanishes for $M=\kappa(\fq)$, where
$\fq\subset B_{(m+n)}$ is any prime ideal. Fix such $\fq$, and
set $\fp:=\fq\cap A_{(m)}$; if $\fp=0$, then $M$ is a $K$-algebra,
where $K$ is the field of fractions of $A_{(m)}$; now, $K$ is a
field of characteristic zero, and $B':=B_{(m+n)}\otimes_{A_{(m)}}K$
is a regular local $K$-algebra, so the induced morphism
$\Spec\,B'\to\Spec\,K$ is regular; since
$H(M)=H_1(\L_{B'/K}\otimes_{B'}M)$, the assertion follows from
corollary \ref{cor_reg-maps-crit}. Notice that this argument
applies especially to the case where $m=0$; for the general case,
we argue by induction on $m$. Hence, suppose that $n\in\N$, $m>0$,
and that the assertion is already known for $\phi_{n,m-1}$.

Consider first the case where $A$ is a discrete valuation ring,
and $\fp$ contains the maximal ideal of $A$, and set
$\bar B:=B\otimes_A\kappa$. Since $\phi_{m,n}$ is flat,
and since $M$ is a $\bar B_{(m+n)}$-module, we have a
natural isomorphism
$$
H(M)\isom
H_1(\L_{\bar B_{(m+n)}/\kappa_{(m)}}\otimes_{\bar B_{(m+n)}}M).
$$
Then the sought vanishing follows from the foregoing, since
$\bar B$ is a formally smooth $\kappa$-algebra (for the preadic
topology of its maximal ideal).

Lastly, suppose that either $A$ is a field, or $\fp$ does not
contain the maximal ideal of $A$ (and $\fp\neq 0$). In either
of these two cases, we may find $f\in\fp$ whose image in
$\kappa_{(m)}$ is not zero. According to
\cite[Ch.VII, n.7, Lemme 3]{BouAC} and
\cite[Ch.VII, n.8, Prop.6]{BouAC}, we may find an automorphism
$\sigma$ of the $A$-algebra $A_{(m)}$ and elements
$g\in A_{(m-1)}[T_m]$, $u\in A_{(m)}^\times$ such that
$\sigma(f)=u\cdot g$, and $g=T_m^d+a_1T_m^{d-1}+\cdots+a_d$
for some $d\geq 0$ and certain elements $a_1,\dots,a_d$ of
the maximal ideal of $A_{(m-1)}$.
Denote by $\sigma':B_{(m)}\isom B_{(m)}$ the $T$-adic completion
of $\sigma\otimes_AB$, and let $\sigma_B:B_{(m+n)}\isom B_{(m+n)}$
be the $T$-adically continuous automorphism that restricts to
$\sigma'$ on $B_{(m)}$, and such that $\sigma_B(T_i)=T_i$ for
$i=m+1,\dots,m+n$. Set $M':=B_{(m+n)}/\sigma_B(\fq)$; by
construction, we have a commutative diagram of $A$-algebras
$$
\xymatrix{
A_{(m)} \ar[rr]^-{\phi_{m,n}} \ar[d]_\sigma & &
B_{(m+n)} \ar[d] \ar[d]^{\sigma_B} \\
A_{(m)} \ar[rr]^-{\phi_{m,n}} & & B_{(m+n)}
}$$
inducing an isomorphism
$$
\L_{B_{(m+n)}/A_{(n)}}\otimes_{B_{(m+n)}}M\isom
\L_{B_{(m+n)}/A_{(n)}}\otimes_{B_{(m+n)}}M'
\qquad
\text{in $\sD(A\Mod)$}.
$$
Thus, we may replace $M$ by $M'$, and assume from start
that $g\in\fp$. In this case, set $B':=B[[T_{m+1},\dots,T_{m+n}]]$,
and notice that both of the natural maps
$$
A_{(m-1)}[T_m]/gA_{(m-1)}[T_m]\to A_{(m)}/gA_{(m)}
\qquad
B'_{(m-1)}[T_m]/gB'_{(m-1)}[T_m]\to B_{(m+n)}/gB_{(m+n)}
$$
are isomorphisms. Since both $\phi_{m,n}$ and the map
$A_{(m-1)}[T_m]\to B'_{(m-1)}[T_m]$ induced by $\phi$ are
flat ring homomorphisms, there follows a natural isomorphism
of $B_{(m+n)}$-modules :
$$
\begin{aligned}
\L_{B_{(m+n)}/A_{(n)}}\otimes_{B_{(m+n)}}M\isom\: &
\L_{B'_{(m-1)}[T_m]/A_{(m-1)}[T_m]}\otimes_{B'_{(m-1)}[T_m]}M \\
\isom\: &
\L_{B'_{(m-1)}/A_{(m-1)}}\otimes_{B'_{(m-1)}}M
\end{aligned}
$$
(\cite[Ch.II, Prop.2.2.1]{Il}).
However, the resulting map $A_{(m-1)}\to B'_{(m-1)}$ is none
else than $\phi_{m-1,n}$, up to a relabeling of the variables;
the vanishing of $H(M)$ then follows from the inductive
assumption (and from corollary \ref{cor_reg-maps-crit}).
\end{pfclaim}

\begin{claim}\label{cl_fin-type-H}
In the situation of the theorem, suppose
furthermore that $A$ and $B$ are complete, and let $M$ be
a $B$-module of finite type. Then $H_1(\L_{B/A}\otimes_BM)$
is a $B$-module of finite type.
\end{claim}
\begin{pfclaim} Let $\bar f:A_0\to\kappa$ be a surjective
ring homomorphism, with $A_0$ a Cohen ring
(\cite[Ch.0, Th.19.8.6(ii)]{EGAIV}), and set
$\bar B:=B\otimes_A\kappa$; in light of
\cite[Ch.0, Lemme 19.7.1.3]{EGAIV}, there exists a flat
local, complete and noetherian $A_0$-algebra $B_0$ fitting
into a cocartesian diagram :
$$
\xymatrix{ A_0 \ar[r]^-\psi \ar[d]_{\bar f} & B_0 \ar[d]^{\bar f_B} \\
          \kappa \ar[r] & \bar B.
}$$
Denote by $\kappa_B$ the residue field of $B$; there follow
natural isomorphisms of $B$-modules :
$$
H_1(\L_{B/A}\otimes_B\kappa_B)\isom
H_1(\L_{\bar B/\kappa}\otimes_{\bar B}\kappa_B)\isom
H_1(\L_{B_0/A_0}\otimes_{B_0}\kappa_B)
$$
(\cite[Ch.II, Prop.2.2.1]{Il})
which, according to proposition \ref{prop_fsmooth-cot}, imply
that $\psi$ is formally smooth (for the preadic topologies
defined by the maximal ideals). By \cite[Ch.0, Th.19.8.6(i)]{EGAIV},
$\bar f$ lifts to a ring homomorphism $f:A_0\to A$, whence a
commutative diagram
$$
\xymatrix{
A_0 \ar[r]^-\psi \ar[d]_{\phi\circ f} & B_0 \ar[d]^{\bar f_B} \\
B \ar[r] & \bar B.
}$$
Then, by \cite[Ch.0, Cor.19.3.11]{EGAIV}, the map $\bar f_B$
lifts to a ring homomorphism $f_B:B_0\to B$. Notice now that
both $f$ and $f_B$ are local maps and induce isomorphisms on
the residue fields; it follows easily that, for suitable $m,n\in\N$,
they extend to surjective maps $g$ and $g_B$ fitting into a
commutative diagram
$$
\xymatrix{
A_{0,(m)} \ar[rr]^-{\psi_{(m,n)}} \ar[d]_g & &
B_{0,(m+n)} \ar[d]^{g_B} \\
A \ar[rr]^-\phi & & B
}$$
whence exact sequences (\cite[Ch.II, Prop.2.1.2]{Il})
$$
\begin{aligned}
& H_1(\L_{B_{0,(m+n)}/A_{0,(m)}}\otimes_{B_{0,(m+n)}}M)\to
H_1(\L_{B/A_{0,(m)}}\otimes_BM)\to
H_1(\L_{B_{0,(m+n)}/B}\otimes_{B_{0,(m+n)}}M) \\
& H_1(\L_{B/A_{0,(m)}}\otimes_BM)\to H_1(\L_{B/A}\otimes_BM)\to
\Omega^1_{A/A_{0,(m)}}\otimes_AM=0.
\end{aligned}
$$
However, claim \ref{cl_m-and-n-reg} applies to $\psi$, and
together with corollary \ref{cor_reg-maps-crit}, it implies
that the first module of the first of these sequences vanishes;
on the other hand, it is easily seen that the third $B$-module
of the same sequence is finitely generated, so the same holds
for the middle term. By inspecting the second exact sequence,
the claim follows.
\end{pfclaim}

We may now conclude the proof of the theorem : let $A^\wedge$
and $B^\wedge$ be the completions of $A$ and $B$; since
$\phi$ is formally smooth, the same holds for its completion
$\phi^\wedge:A^\wedge\to B^\wedge$. First, we show that
$\Spec\,\phi^\wedge$ is regular; to this aim, it suffices
to check that
$H^\wedge(M):=H_1(\L_{B^\wedge/A^\wedge}\otimes_{B^\wedge}M)$
vanishes for every $B^\wedge$-module $M$ (corollary
\ref{cor_reg-maps-crit}). However, we know already that
$H^\wedge(M)$ is a $B^\wedge$-module of finite type, if
the same holds for $M$ (claim \ref{cl_fin-type-H}); especially,
for such $M$, $H^\wedge(M)$ is separated for the adic topology
of $B^\wedge$ defined by the maximal ideal. Moreover,
$H^\wedge(M)=0$, if $M$ is the residue field of $B^\wedge$
(proposition \ref{prop_fsmooth-cot}). Then the assertion follows
from claim \ref{cl_general-remark}. Lastly, the natural morphism
$\Spec\,A^\wedge\to\Spec\,A$ is regular by assumption, hence
the same holds for the induced morphism $\Spec\,B^\wedge\to\Spec\,A$
(lemma \ref{lem_yoga-reg}(i)). Finally, since $B^\wedge$ is
a faithfully flat $B$-algebra, we conclude that $\Spec\,\phi$
is regular, by virtue of lemma \ref{lem_yoga-reg}(ii).
\end{proof}

\begin{proposition}\label{prop_block-buster}
Let $A$ be a noetherian local ring, and $\phi:A\to B$
an ind-{\'e}tale local ring homomorphism. Then :
\begin{enumerate}
\item
$A$ is quasi-excellent if and only if the same holds for $B$.
\item
If $A$ is excellent, the same holds for $B$.
\end{enumerate}
\end{proposition}
\begin{proof} Set $f:=\Spec\,\phi$, and
$f^\wedge:=\Spec\,\phi^\wedge$, where
$\phi^\wedge:A^\wedge\to B^\wedge$ is the map of complete
local rings obtained from $\phi$. Notice first that, for
every $y\in\Spec\,A$, and every $x\in f^{-1}(y)$, the stalk
$\cO_{\!f^{-1}y,x}$ is an ind-{\'e}tale $\kappa(y)$-algebra,
{\em i.e.} is a separable algebraic extension of $\kappa(y)$,
hence $f^{-1}(y)$ is geometrically regular, and therefore $f$
is regular and faithfully flat. Moreover, $f$ is the limit of
a cofiltered system of formally \'etale morphisms, hence it
is formally \'etale; {\em a fortiori}, $B$ is also formally
smooth over $A$ for the preadic topologies, so $B^\wedge$ is
formally smooth over $A^\wedge$ for the preadic topologies,
and consequently $f^\wedge$ is a regular morphism (theorem
\ref{th_block-buster}).

(i): If $B$ is quasi-excellent, claim \ref{cl_descent-G-ring}
and lemma \ref{lem_Nagatas}(iv) imply that $A$ is quasi-excellent.
Next, assume that $A$ is quasi-excellent; we have a commutative
diagram \eqref{eq_f-and-f-wedge}, in which $\pi_A$ is a regular
morphism, and then the same holds for
$\pi_A\circ f^\wedge=f\circ\pi_B$ (lemma \ref{lem_yoga-reg}(i)).
Now, let $y\in\Spec\,B$ be any point, and set $z:=f(y)$; we
know that $(f\circ\pi_B)^{-1}(z)$ is a geometrically regular
affine scheme, so write it as the spectrum of a noetherian
ring $C$. The stalk $E:=\cO_{\!f^{-1}(z),y}$ is a separable
algebraic field extension of $\kappa(z)$, hence
$\pi_B^{ -1}(y)\simeq\Spec\,C\otimes_BE$, the spectrum of a
localization of $C$, so it is again geometrically regular,
{\em i.e.} $\pi_B$ is a regular morphism, and we conclude
by lemma \ref{lem_Nagatas}(iii,iv).

(ii): By (i), it remains only to show that $B$ is
universally catenarian, provided the same holds for $A$.
To this aim, it suffices to apply \cite[Ch.IV, Lemma 18.7.5.1]{EGA4}.
\end{proof}

\begin{proposition}\label{prop_jump-of-the-way}
In the situation of \eqref{subsec_again-Frobby}, suppose that
$A$ is noetherian and $\Phi_A$ is a finite ring homomorphism,
and for every prime ideal $\fp\subset A$, set
$\kappa(\fp):=A_\fp/\fp A_\fp$ (the residue field of the point
$\fp\in\Spec\,A$). Then
$$
\dim A_\fp/\fq A_\fp=\dim_{\kappa(\fp)}\Omega^1_{\kappa(\fp)/\F_p}-
\dim_{\kappa(\fq)}\Omega^1_{\kappa(\fq)/\F_p}
$$
for every pair of prime ideals $\fq\subset\fp\subset A$.
\end{proposition}
\begin{proof} To begin with, we notice :

\begin{claim}\label{cl_stability}
In the situation of \eqref{subsec_again-Frobby},
suppose that $\Phi_A$ is a finite map. We have :
\begin{enumerate}
\item
For every $A$-algebra $B$ of essentially finite type,
$\Phi_B$ is finite as well.
\item
Suppose moreover that $A$ is noetherian, $I\subset A$
is any ideal, and let $A^\wedge_I$ be the $I$-adic
completion of $A$. Then :
\begin{enumerate}
\item
$\Phi_{A_I^\wedge}=\one_{A_I^\wedge}\otimes_A\Phi_A$ is finite,
and $\Omega^1_{A^\wedge_I/\F_p}=A^\wedge_I\otimes_A\Omega^1_{A/\F_p}$.
\item
If $A$ is reduced, the same holds for $A^\wedge_I$.
\end{enumerate}
\item
Suppose that $(A,\fm_A)$ is a local noetherian domain, and
denote by $K$ (resp. $\kappa$) the field of fractions
(resp. the residue field) of $A$. Then
\set\begin{equation}\label{eq_take-diff}
\dim_K\Omega^1_{K/\F_p}=\dim_\kappa\Omega^1_{\kappa/\F_p}+\dim A.
\end{equation}
\end{enumerate}
\end{claim}
\begin{pfclaim}(i): Suppose first that $B=A[X]$; then it is
easily seen that $\Phi_B=\Phi_A\otimes_{\F_p}\Phi_{\F_p[X]}$,
whence the contention, in this case.
By an easy induction, we deduce that the claim holds as well
for any free polynomial $A$-algebra of finite type. Next, let
$I\subset A$ be any ideal; since $\Phi_A(I)\subset I$, it is
easily seen that $\Phi_{A/I}=\Phi_A\otimes_AA/I$, so $\Phi_{A/I}$
is finite, and therefore the assertion holds for any $A$-algebra
of finite type. Lastly, let $S\subset A$ be any multiplicative
subset; since $\Phi_A(S)\subset S$, it is easily seen that
\set\begin{equation}\label{eq_loc_Kunz}
\Phi_{S^{-1}A}=S^{-1}\Phi_A
\end{equation}
and the assertion follows.

(ii.a): Set $B:=A_{(\Phi)}$; by assumption, $\Phi_A:A\to B$
is a finite $A$-linear map, hence its $I$-adic completion
$(\Phi_A)^\wedge_I:A_I^\wedge\to B_I^\wedge$ equals
$\one_{A^\wedge_I}\otimes_A\Phi_A$ (\cite[Th.8.7]{Mat}).
Say that $I$ is generated by $r$ elements of $A$; then it
is easily seen that $I^{(p-1)r+1}\subset\Phi_A(I)\subset I$,
so we have a natural $A^\wedge_I$-linear identification on
$I$-adic completions :
\set\begin{equation}\label{eq_corinna}
B^\wedge_I\isom(A^\wedge_I)_{(\Phi)}
\end{equation}
and under this identification, $(\Phi_A)^\wedge_I=\Phi_{A^\wedge_I}$,
whence the first assertion of (ii.a). Next, we notice that
\eqref{eq_corinna} yields a natural identification
$\Omega^1_{A^\wedge_I/\F_p}=\Omega^1_{B^\wedge_I/\F_p}$, and on the
other hand, the sequence of ring homomorphisms
$$
\F_p\to A^\wedge_I\xrightarrow{\ \Phi_{A^\wedge_I}\ } B_I^\wedge
$$
yields natural identifications :
$$
\Omega^1_{B_I^\wedge/\F_p}=
\Omega^1_{B_I^\wedge/A_I^\wedge}=
A_I^\wedge\otimes_A\Omega^1_{B/A}=
A^\wedge_I\otimes_A\Omega^1_{B/\F_p}=
B_I^\wedge\otimes_B\Omega^1_{B/\F_p}=
A^\wedge_I\otimes_A\Omega^1_{A/\F_p}
$$
where the last equality is induced by \eqref{eq_corinna} and
the identity map $A\isom B$. This completes the proof of the
second assertion.

(ii.b): Since $A$ is reduced, $\Phi_A$ is injective, and then
the same holds for its completion $(\Phi_A)_I^\wedge$. But we
have seen that the latter is naturally identified with
$\Phi_{A^\wedge_I}$, whence the claim.

(iii): Notice that $\Omega^1_{K/\F_p}$ is a $K$-vector space
of finite dimension, since its dimension equals $[K:K^p]$
(\cite[Ch.0, Th.21.4.5]{EGAIV}),
and the latter is finite, by (i); the same argument applies
to the $\kappa$-vector space $\Omega^1_{\kappa/\F_p}$. Denote
by $\delta(A)$ the difference between the left and right
hand-side of \eqref{eq_take-diff}; we have to show that
$\delta(A)=0$. Let $A^\wedge$ be the $\fm_A$-adic completion
of $A$, and $\fp\subset A^\wedge$ any minimal
prime ideal such that $\dim A^\wedge/\fp=d:=\dim A$. In view
of (ii.b), $L:=(A^\wedge)_\fp$ is a field; taking into account
(ii.a), we get
$$
\Omega^1_{L/\F_p}=L\otimes_{A^\wedge}\Omega^1_{A^\wedge/\F_p}=
L\otimes_A\Omega^1_{A/\F_p}=L\otimes_K\Omega^1_{K/\F_p}.
$$
Hence $\dim_L\Omega^1_{L/\F_p}=\dim_K\Omega^1_{K/\F_p}$, so
we see that
\set\begin{equation}\label{eq_delta-fp}
\delta(A)=\delta(A^\wedge/\fp)
\qquad
\text{for any prime ideal $\fp\subset A^\wedge$
such that $\dim A^\wedge/\fp=d$}.
\end{equation}
We may then replace $A$ by $A^\wedge/\fp$, after which we
may assume that $A$ is a complete local domain. In this
case, $A$ contains a field mapping isomorphically onto
$\kappa$, and there is a finite map of $\kappa$-algebras
$A':=\kappa[[T_1,\dots,T_d]]\to A$ (\cite[Th.29.4(iii)]{Mat}).
Denote by $K'$ the field of fractions of $A'$, and notice that
$[K:K']=[K:K^p]\cdot[K^p:K']=[K:K^p]\cdot[K^p:K'{}^p]$;
on the other hand, $\Phi_K$ induces an isomorphism $K\isom K^p$,
hence have as well $[K:K']=[K^p:K'{}^p]$, so $[K:K^p]=[K':K'{}^p]$
and finally :
$$
\dim_{K'}\Omega^1_{K'/\F_p}=\dim_K\Omega^1_{K/\F_p}
$$
(\cite[Ch.0, Cor.21.2.5]{EGAIV}). Since $\dim A'=\dim A$, we
conclude that $\delta(A)=\delta(A')$.
Hence, it suffices to check that $\delta(A')=0$. Furthermore,
set $B:=\kappa[T_1,\dots,T_d]$, and let $\fm\subset B$ be the
maximal ideal generated by $T_1,\dots,T_d$; we have
$A'=B^\wedge_\fm$, so \eqref{eq_delta-fp} further reduces
to showing that $\delta(B_\fm)=0$, which shall be left as
an exercise for the reader.
\end{pfclaim}

Now, in view of claim \ref{cl_stability}(i), we may replace
$A$ by $A_\fp/\fq A_\fp$, and assume from start that $A$ is
a local domain, that $\fq=0$ and that $\fp$ is the maximal
ideal; in this case, the sought identity is given by
claim \ref{cl_stability}(iii).
\end{proof}

\begin{theorem}\label{th_Kunz-exc}
With the notation of \eqref{subsec_again-Frobby}, suppose that
$A$ is noetherian. We have :
\begin{enumerate}
\item
If\/ $\Phi_{\!A}$ is a finite ring homomorphism, then $A$ is
excellent.
\item
Conversely, suppose that $A$ is a local Nagata ring, with
residue field $k$, and that $[k:k^p]$ is finite. Then
$\Phi_{\!A}$ is a finite ring homomorphism.
\end{enumerate}
\end{theorem}
\begin{proof} (i): We reproduce the proof from
\cite[Appendix, Th.108]{Mat-old}.

\begin{claim}\label{cl_Kunz-quasi-exc}
If $\Phi_{\!A}$ is finite, then $A$ is quasi-excellent.
\end{claim}
\begin{pfclaim} We check first the openness condition
for the regular loci. To this aim, in light of claim
\ref{cl_stability}(i), we may assume that $A$ is an integral
domain, and it suffices to prove that the regular locus of
$\Spec\,A$ is an open subset.
Set $B:=A_{(\Phi)}$ (notation of \eqref{subsec_again-Frobby}),
and let $\fp\subset A$ be any prime ideal; by  theorem
\ref{th_Kunz-by-Matsu} and \eqref{eq_loc_Kunz}, the ring
$A_\fp$ is regular if and only if $(\Phi_A)_\fp$ is a flat
ring homomorphism, if and only if $B_\fp$ is a flat
$A_\fp$-module. Since, by assumption, $B$ is a finite
$A$-module, the latter holds if and only if $B_\fp$ is
a free $A_\fp$-module. Then, the assertion follows from
\cite[Th.4.10]{Mat}.

Next, we show that $A$ is a G-ring. By claim \ref{cl_stability}(i),
we may assume that $A$ is local, and we let $A^\wedge$
be the completion of $A$. From claim
\ref{cl_stability}(ii.a) we see that the natural map
$$
A_{(\Phi)}\derotimes_AA^\wedge\to(A^\wedge)_{(\Phi)}
$$
is an isomorphism in $\sD(A\Mod)$. Then the assertion follows
from \cite[Lemma 6.5.13(i)]{Ga-Ra} and corollary \ref{cor_crit-q.excel}.
\end{pfclaim}

Now, it follows easily from proposition \ref{prop_jump-of-the-way}
and claim \ref{cl_stability}(i), that $A$ is universally catenarian;
combining with claim \ref{cl_Kunz-quasi-exc}, we obtain (i).

(ii): The assertion to prove is that $A_{(\Phi)}$ is a finite
$A$-module. To this aim, denote by $I\subset A$ the nilradical
ideal; we remark :

\begin{claim}\label{cl_reduce-to-red}
It suffices to show that $(A/I)_{(\Phi)}$ is a finite
$A$-module.
\end{claim}
\begin{pfclaim} Indeed, suppose that $(A/I)_{(\Phi)}$ is a
finite $A$-module; we shall deduce, by induction on $s$,
that $(A/I^s)_{(\Phi)}$ is a finite $A$-module, for every
$s>0$. Since $A$ is noetherian, we have $I^t=0$ for $t>0$
large enough, so the claim will follow. The assertion for
$s=1$ is our assumption. Suppose therefore that $s>1$, and
that we already know the assertion for $s-1$. Notice that
$I^s/I^{s-1}$ is a finite $A/I$-module, hence a finite
$(A/I)_{(\Phi)}$-module, by our assumption. Since
$(A/I^{s-1})_{(\Phi)}$ is a finite $A$-module, we deduce
easily that the assertion holds for $s$, as required.
\end{pfclaim}

In view of claim \ref{cl_reduce-to-red}, we may replace
$A$ by $A/I$ (which is obviously still a Nagata ring), and
assume from start that $A$ is reduced. In this case, let
$\{\fp_1,\dots,\fp_t\}$ be the set of minimal prime ideals
of $A$; we have a commutative diagram of ring homomorphisms :
$$
\xymatrix{ A \ar[rrr]^-{\Phi_A} \ar[d] & & & A \ar[d] \\
\prod_{i=1}^tA/\fp_i \ar[rrr]^-{\prod_{i=1}^t\Phi_{A/\fp_i}}
& & & \prod_{i=1}^tA/\fp_i
}$$
whose vertical arrows are injective and finite maps. Suppose
now that $(A/\fp_i)_{(\Phi)}$ is a finite $A/\fp_i$-module
for every $i=1,\dots,t$. Then $\prod_{i=1}^t(A/\fp_i)_{(\Phi)}$
is a finite $A$-module, and therefore the same holds for its
$A$-submodule $A_{(\Phi)}$. Thus, we are reduced to checking
the sought assertion for each of the quotients $A/\fp_i$ (which
are still Nagata rings), and hence we may assume from start that
$A$ is a domain. In this case, denote by $K$ the field of
fractions of $A$; by the definition of Nagata ring, we are
further reduced to checking that $\Phi_K:K\to K$ is a finite
field extension, {\em i.e.} that $[K:K^p]$ is finite.
Denote $A^\wedge$ the completion of $A$; we remark :

\begin{claim}\label{cl_Zariski-Nagata}
(i)\ \ The endomorphism $\Phi_{\!A^\wedge}$ is a finite map.
\begin{enumerate}
\addenu
\item
$\Spec\,K\otimes_AA^\wedge$ is a geometrically reduced
$K$-scheme.
\end{enumerate}
\end{claim}
\begin{pfclaim}(i): Indeed, $A^\wedge$ is a quotient of a
power series ring $B:=k[[T_1,\dots,T_r]]$
(\cite[Ch.0, Th.19.8.8(i)]{EGAIV}), hence it
suffices to show that $\Phi_B$ is finite. However,
$B^p=k^p[[T_1^p,\dots,T_r^p]]$, and $[k:k^p]$ is finite
by assumption, whence the claim.

(ii): This is \cite[Ch.IV, Th.7.6.4]{EGAIV-2}.
\end{pfclaim}

Let $\fp\subset A^\wedge$ be any minimal prime ideal;
it follows from claim \ref{cl_Zariski-Nagata}(ii) that
$L:=(A^\wedge)_\fp$ is a separable field extension of $K$
(\cite[Ch.IV, Prop.4.6.1]{EGAIV-2}), so the natural map
$$
\Omega^1_{K/\F_p}\otimes_KL\to\Omega^1_{L/\Z}
$$
is injective (\cite[Ch.0, Cor.20.6.19(i)]{EGAIV}); on the
other hand, claim \ref{cl_Zariski-Nagata}(i) implies that
$[L:L^p]$ is finite, {\em i.e.} $\Omega^1_{L/\F_p}$ is
a finite-dimensional $L$-vector space
(\cite[Ch.0, Th.21.4.5]{EGAIV}). We conclude that
$\Omega^1_{K/\F_p}$ is a finite-dimensional $K$-vector
space, whence the contention, again by {\em loc.cit.}
\end{proof}

\subsection{Normalization, weak normalization and
\texorpdfstring{$p$}{p}-integral closure}
Recall that for an inclusion of rings $A\subset B$, we
denote by $\mathrm{i.c.}(A,B)$ the integral closure of
$A$ into $B$ (see definition \ref{def_ic}(i)); recall
moreover that
\set\begin{equation}\label{eq_normal-and-localization}
\mathrm{i.c.}(S^{-1}A,S^{-1}B)=S^{-1}\mathrm{i.c.}(A,B)
\qquad
\text{for every subset $S\subset A$}.
\end{equation}
More generally, let $X$ be a scheme, $\cA\subset\cB$ an
inclusion of quasi-coherent $\cO_X$-algebras; then the rule :
$U\mapsto\cC(U):=\mathrm{i.c.}(\cA(U),\cB(U))$ for every
affine open subset $U\subset X$ defines a subpresheaf
$\cC\subset\cB$ on the site of affine open subsets of $X$,
and \eqref{eq_normal-and-localization} easily implies that
$\cC$ is a sheaf on this site, so it extend uniquely on the
Zariski site of $X$ to a quasi-coherent $\cO_X$-subalgebra
of $\cB$, which we call the {\em integral closure of $\cA$
in $\cB$}, and denote
$$
\mathrm{i.c.}(\cA,\cB).
$$
Also, for every ring $A$, we let $\cN_A$ be the nilradical
of $A$, and set $A_\red:=A/\cN_A$, the maximal reduced quotient
of $A$. Every ring homomorphism $f:A\to A'$ induces a ring
homomorphism $f_\red:A_\red\to A'_\red$. Recall that
$$
\cN_{S^{-1}A}=S^{-1}A
\qquad\text{and}\qquad
S^{-1}(A_\red)=(S^{-1}A)_\red
\qquad
\text{For every subset $S\subset A$}.
$$
Arguing as in the foregoing, we deduce that for every
quasi-coherent $\cO_X$-algebra $\cA$, the rule :
$U\mapsto\cN_{\cA(U)}$ for every affine open subset
$U\subset X$ extends uniquely to a quasi-coherent ideal
$$
\cN_\cA\subset\cA
$$
called the {\em nilradical of $\cA$}. The $\cO_X$-algebra
$\cA_\red\!:=\!\cA\!/\!\cN_\cA$ is quasi-coherent, and every
morphism $\phi:\cA\to\cA'$ of quasi-coherent $\cO_X$-algebras
induces a morphism $\phi_\red:\cA_\red\to\cA'_\red$.

\begin{lemma}\label{lem_Krullu}
Let $A$ be a Krull domain, $F$ the field of fractions
of $A$, $B$ a flat $A$-algebra, and suppose that
$B\otimes_A\kappa(\fp)$ is reduced, for every prime ideal
$\fp\in\Spec\,A$ of height one. Then $B$ is integrally
closed in $B\otimes_AF$.
\end{lemma}
\begin{proof} See \cite[\S12]{Mat} for the basic generalities
on Krull domains; especially, \cite[Th.12.6]{Mat} asserts
that if $A$ is a Krull domain, the natural sequence of $A$-modules :
$$
\cE\quad :\quad
0\to A\to F\to\bigoplus_{\mathrm{ht}\,\fp=1}F/A_\fp\to 0
$$
is short exact, where $\fp\in\Spec\,A$ ranges over the
prime ideals of height one. By flatness, $\cE\otimes_AB$
is still exact, hence $B=\bigcap_{\mathrm{ht}\,\fp=1}B\otimes_AA_\fp$
(where the intersection takes place in $B\otimes_AF$).
We are therefore reduced to the case where $A$ is a discrete
valuation ring. Let $t$ denote a chosen generator of the
maximal ideal of $A$, and suppose that $x\in B\otimes_AF$
is integral over $B$, so that $x^n+b_1x^{n-1}+\cdots+b_n=0$
in $B\otimes_AF$, for some $b_1,\dots,b_n\in B$; let also
$r\in\N$ be the minimal integer such that we have $x=t^{-r}b$
for some $b\in B$. We have to show that $r=0$; to this aim,
notice that $b^n+t^rb_1b^{n-1}+\cdots+t^{rn}b_n=0$ in $B$;
if $r>0$, it follows that the image $\bar b$ of $b$ in $B/tB$
satisfies the identity : $\bar b{}^n=0$, therefore $\bar b=0$,
since by assumption, $B/tB$ is reduced. Thus $b=tb'$ for some
$b'\in B$, and $x=t^{1-r}b'$, contradicting the minimality of $r$.
\end{proof}

The following proposition generalizes
\cite[Ch.IV, Prop.6.14.4]{EGAIV-2}.

\begin{proposition}\label{prop_normal-and-geom-red-bc}
Let $f:X\to S$ be a flat morphism of schemes with geometrically
reduced fibres, $\cB_1\subset\cB_2$ two quasi-coherent
$\cO_S$-algebras, and $\cB_3:=\mathrm{i.c.}(\cB_1,\cB_2)$.
Suppose that :
\begin{enumerate}
\alphaenu
\item
either $f$ is locally finitely presented,
\item
or else, $S$ is locally noetherian.
\end{enumerate}
\item
Then $f^*\cB_3=\mathrm{i.c.}(f^*\cB_1,f^*\cB_2)$ and
$f^*(\cB_1)_\red=f^*(\cB_{1,\red})$.
\end{proposition}
\begin{proof}The assertion is local on both $S$ and $X$,
hence we may assume that $S=\Spec\,A$ and $X=\Spec\,A'$,
for a ring $A$ and an $A$-algebra $A'$. Then also $\cB_1$
and $\cB_2$ are the quasi-coherent algebras associated
with $A$-algebras $B_1\subset B_2$, and $\cB_3$ is the
quasi-coherent algebra associated with
$B_3:=\mathrm{i.c.}(B_1,B_2)$.
Suppose first that condition (a) holds, and write $A$ as the
colimit of a filtered system $(A_\lambda~|~\lambda\in\Lambda)$
of noetherian subrings; we remark :

\begin{claim}\label{cl_descend-geom-reduced}
Let $A_\bullet:=(A_\lambda~|~\lambda\in\Lambda)$ be a filtered
system of rings, $A$ the colimit of $A_\bullet$, and
$X\to S:=\Spec\,A$ a flat morphism of schemes of finite
presentation, with geometrically reduced fibres. Then for some
$\lambda\in\Lambda$ there exists a flat morphism of schemes
with geometrically reduced fibres $X_\lambda\to S_\lambda$, and
an isomorphism of $S$-schemes $S\times_{S_\lambda}X_\lambda\isom X$.
\end{claim}
\begin{pfclaim} By \cite[Ch.IV, Th.8.8.2(ii)]{EGAIV-3}, we may
find $\lambda\in\Lambda$, a finitely presented morphism
$f_\lambda:X_\lambda\to S_\lambda$, and an isomorphism of
$S$-schemes : $X\isom X_\lambda\times_{S_\lambda}S$. We set
$S_\mu:=\Spec\,A_\mu,X_\mu:=X_\lambda\times_{S_\lambda}S_\mu$, and
$f_\mu:=S_\mu\times_{S_\lambda}f_\lambda:X_\mu\to S_\mu$ for every
$\mu\in\Lambda$ with $\mu\geq\lambda$; after replacing
$\Lambda$ by a cofinal subset, we may then assume that
$f_\mu$ is defined for every $\mu\in\Lambda$, and moreover that
each such $f_\mu$ is flat (\cite[Ch.IV, Cor.11.2.6.1]{EGAIV-3}).
For every $\mu\in\Lambda$, let $Z_\mu\subset S_\mu$ be the subset
of all $s\in S_\mu$ such that the fibre $f_\mu^{-1}(s)$ is not
geometrically reduced; by \cite[Ch.IV, Th.9.7.7]{EGAIV-3}, $Z_\mu$
is a constructible subset of $S_\mu$. By assumption, we have
$$
\bigcap_{\mu\in\Lambda}p^{-1}_\mu Z_\mu=\emptyset
$$
and it is clear that $p^{-1}_\mu Z_\mu\subset p^{-1}_\lambda Z_\lambda$
whenever $\mu\geq\lambda$. Then, $Z_\mu=\emptyset$ for some
$\mu\in\Lambda$ (\cite[Ch.IV, Prop.1.8.2, Cor.1.9.8]{EGAIV}),
hence we may replace $\Lambda$ by a still smaller cofinal
subset, and achieve that all the $f_\lambda$ have geometrically
reduced fibres.
\end{pfclaim}

By claim \ref{cl_descend-geom-reduced}, we find
$\lambda\in\Lambda$ and a flat ring homomorphism
$\phi_0:§A_{\lambda_0}\to A'_0$ such that :
\begin{itemize}
\item
all the fibres of $\Spec(\phi_0)$ are geometrically reduced
\item
there exists an isomorphism of $A$-algebras
$A'\isom A\otimes_{A_\lambda}A'_0$.
\end{itemize}
We need to show that
$A'\otimes_AB_3=\mathrm{i.c.}(A'\otimes_AB_1,A'\otimes_AB_2)$
and $A'\otimes_AB_{1,\red}=(A'\otimes_AB_1)_\red$. However, in
view of the natural identifications :
$A'\otimes_AB_i\isom A'_0\otimes_{A_\lambda}B_i$ for $i=1,2,3$, we
may replace $A$ and $A'$ by $A_\lambda$ and $A'_0$, and assume
that $A$ is noetherian, {\em i.e.} that condition (b) holds.
For the first identity, arguing as in the proof of
\cite[Ch.IV, Prop.6.14.4 and Cor.6.14.5]{EGAIV-2},
we reduce to the following two separate cases :
\begin{itemize}
\item[(a')]
$B_1$ is a domain of finite type over $A$ and $B_2$ is
a finite extension of $\Frac(B_1)$
\item[(b')]
$B_1$ and $B_2$ are fields, and $B_1$ is algebraically
closed in $B_2$ (so $B_3=B_1$ in this case).
\end{itemize}
In case (a') holds, $B_3$ is a Krull domain (\cite[Th.33.10]{Na}),
and the assertion follows from lemma \ref{lem_Krullu}. In case
(b'), we may replace $A$ by $B_1$ and $A'$ by $A'\otimes_AB_1$,
and assume that $A$ is a field, $A'$ is a geometrically
reduced $A$-algebra, $B:=B_2$ is a field extension of $A$ such
that $A$ is algebraically closed in $B$, and we need to check
that $A'$ is integrally closed in $A'\otimes_AB$. Now, pick
a basis $(b_i~|~i\in I)$ for the $A$-vector space $B$,
with $b_{i_0}=1$ for some $i_0\in I$, and say that
$\sum_{i\in J}a_i\otimes b_i\in A'\otimes_AB$ is integral
over $A'$ (for some finite subset $J\subset I$); we
need to check that $a_i=0$ for every $i\in J\setminus\{i_0\}$.
However, since $A'$ is reduced, it suffices to check that
the image of $a_i$ vanishes in $\kappa(\fp)$, for every
minimal prime ideal $\fp$ of $A'$. Thus, we may further
assume that $A'$ is a field, in which case it is a separable
(not necessarily algebraic) extension of $A$, since it is
a geometrically reduced $A$-algebra. Then, $A'\otimes_AB$ is
a domain, according to \cite[Ch.V, \S17, n.2, Cor.]{Bourbaki},
and $A'$ is algebraically closed in $\Frac(A'\otimes_AB)$, by
\cite[Ch.V, \S17, Exerc.2(b)]{Bourbaki}.

For the second identity, notice that
$A'\otimes_A\cN_{B_1}\subset\cN_{A'\otimes_AB_1}$; hence it suffices
to check that if $B_1$ is reduced, the same holds for
$A'\otimes_AB_1$. Moreover, we may assume that $B_1$ is an
$A$-algebra of finite type, in which case $B_1$ is noetherian,
and it has finitely many minimal prime ideals $\fq_1,\dots,\fq_n$,
and since $B_1$ is reduced, the natural map
$B_1\to\prod_{i=1}^nB_{\fq_i}$ is injective. Then, the same holds
for the induced map
$A'\otimes_AB_1\to\prod_{i=1}^nA'\otimes_AB_{\fq_i}$. However,
$B_{\fq_i}$ is a field, and $A'\otimes_AB_{\fq_i}$ is reduced
for every $i=1,\dots,n$, since $A'$ has geometrically reduced
fibres over $A$; the assertion follows.
\end{proof}

\begin{definition}\label{def_pro-smooth}
Let $f:X'\to X$ be a morphism of schemes, $\xi'$ a geometric
point of $X'$ localized at $x'\in X'$, and set $\xi:=f(\xi')$.
We shall say that $f$ is {\em pro-smooth at the point $x'$\/}
if the induced map of strictly henselian local rings
$$
\cO_{\!X(\xi),\xi}\to\cO_{\!X'(\xi'),\xi'}
$$
is ind-smooth ({\em i.e.} a filtered colimit of smooth
ring homomorphisms).
\end{definition}

\begin{corollary}\label{cor_pro-smooth}
Let $f:X'\to X$ be a morphism of schemes, $\xi'$ a geometric
point of $X'$, $\cR\subset\cS$ a monomorphism of quasi-coherent
$\cO_{\!X}$-algebras, and
$$
\phi:f^*\mathrm{i.c.}(\cR,\cS)\to\mathrm{i.c.}(f^*\cR,f^*\cS)
$$
the induced morphism. Suppose that $f$ is pro-smooth at the
support $x'$ of\/ $\xi'$; then we have :
\begin{enumerate}
\item
$f$ is flat at the point $x'$.
\item
$\phi_{\xi'}$ is an isomorphism.
\end{enumerate}
\end{corollary}
\begin{proof}(i): Set $x:=f(x'),\xi:=f(\xi')$. We have a commutative
diagram of local $S$-schemes
$$
\xymatrix{
X'(\xi') \ar[rr]^-{f_{(\xi')}} \ar[d]_{\psi'} & &
X(\xi) \ar[d]^\psi \\
X'(x') \ar[rr]^-{f_{(x')}} & & X(x)
}$$
whose vertical arrows are ind-\'etale morphisms, and whose top
horizontal arrow is ind-smooth by assumption. Especially, since
$\psi'$ is faithfully flat, the same holds for $f_{(x')}$,
{\em i.e.} $f$ is flat at $x'$. 

(ii): We may assume that $X=\Spec\,A$ for some ring $A$. Let
$g:X'(\xi')\to X'$ be the natural morphism. It suffices to
show that $g^*\phi$ is an isomorphism. However,
$g^*\mathrm{i.c.}(f^*\cR,f^*\cS)=\mathrm{i.c.}(g^*f^*\cR,g^*f^*\cS)$
by proposition \ref{prop_normal-and-geom-red-bc}, and $g^*\phi$ is
the induced morphism
$$
(f\circ g)^*\mathrm{i.c.}(\cR,\cS)\to
\mathrm{i.c.}((f\circ g)^*\cR,(f\circ g)^*\cS).
$$
Notice that $f\circ g$ is induced by an ind-smooth ring
homomorphism $A\to\cO_{\!X',\xi'}$, so we may invoke again
proposition \ref{prop_normal-and-geom-red-bc} to conclude.
\end{proof}

\sset\subsubsection{}\label{subsec_notation}
Recall that a morphism of schemes $\phi:X\to Y$ is called
{\em universally injective} or {\em radicial} if for every
field $K$, the induced map on $K$-valued points $X(K)\to Y(K)$
is injective (\cite[Ch.I, Def.3.5.4]{EGAI}). It is easily
seen that the latter holds if and only if $\phi$ is injective
and the residue field extension $\kappa(x)\to\kappa(\phi(x))$
is purely inseparable for every $x\in X$.
In the rest of this section we study certain classes of ring
homomorphisms $f:A\to B$ such that $\Spec\,f$ is radicial, in
which case we just say that $f$ is {\em radicial}. We shall
see that $f$ is radicial if and only if the kernel $I_f$ of
the multiplication law $B\otimes_AB\to B$ lies in the nilradical
$\cN_{B\otimes_AB}$ (corollary \ref{cor_yanagihara}(i)).

\begin{lemma}\label{lem_trivialities}
{\em(i)}\ \
Let $S$ be a scheme, and $\phi:Y\to X$ and $\phi':Y'\to X'$ two
surjective morphisms of $S$-schemes. Then
$\phi\times_S\phi':Y\times_SY'\to X\times_SX'$ is surjective.

{\em(ii)}\ \
Let $f\!:\!A\!\to\!B$ be a ring homomorphism such that $\Spec\,f$ is
surjective; then $\Ker\,f\subset\cN_A$.

{\em(iii)}\ \
Let $A$ be a ring, $f:B\to C$ and $f':B'\to C'$ two integral
homomorphisms of $A$-algebras with $\Ker\,f\subset\cN_B$ and
$\Ker\,f'\subset\cN_{B'}$. Then
$f\otimes_Af':B\otimes_AB'\to C\otimes_AC'$ is integral and
its kernel lies in $\cN_{B\otimes_AB'}$.
\end{lemma}
\begin{proof}(i): (See \cite[Ch.I, Prop.3.5.2(i)]{EGAI}). We have
$\phi\times_S\phi'=(\phi\times_SX')\circ(Y\times_S\phi')$, so it
suffices to show that both $\phi\times_SX'$ and $Y\times_S\phi'$
are surjective. But notice that the natural isomorphism
$Y\times_SX'\isom Y\times_X(X\times_SX')$ identifies
$\phi\times_SX'$ with $\phi\times_X(X\times_SX')$, and similarly
for $Y\times_S\phi'$; hence, we are reduced to checking if
$\phi:Y\to X$ is a surjective morphism of schemes, and
$\psi:Z\to X$ is any $X$-scheme, then the base change
$\phi\times_XZ:Y\times_XZ\to Z$ is surjective. Thus, let
$y\in Y$ and $z\in Z$ with $x:=\phi(y)=\psi(z)$, and let
$\kappa(z)\leftarrow\kappa(x)\to\kappa(y)$ be the induced
residue field extensions; it suffices to check that
$\Spec\,\kappa(y)\times_{\Spec\,\kappa(x)}\Spec\,\kappa(z)\neq\emptyset$,
which is clear, since $\kappa(y)\otimes_{\kappa(x)}\kappa(z)\neq 0$.

(ii): Let $a\in\Ker\,f$, and $\fp\in\Spec\,A$; by assumption
there exists $\fq\in\Spec\,B$ with $\fp=f^{-1}\fq$, so that
$a\in\fp$, whence the assertion.

(iii): Obviously $f\otimes_Af'$ is integral. Since
$(B\otimes_AB')_\red=(B_\red\otimes_{A_\red}B'_\red)_\red$, and
similarly for $(C\otimes_AC')_\red$, we may replace $f$ and
$f'$ by $f_\red$ and $f'_\red$, and assume from start that $f$ and
$f'$ are injective. Then, by the going up theorem, both $\Spec\,f$
and $\Spec\,f'$ are surjective, hence the same holds for
$\Spec\,f\otimes_Af'$, by (i), and then the assertion follows
from (ii).
\end{proof}

The following proposition \ref{prop_from_SGA6} is borrowed
from \cite[Exp.XII, Lemma 2.6]{SGA6}, and its corollary
\ref{cor_yanagihara}(ii) can also be easily deduced from
\cite[Th.1]{Yan}.

\begin{proposition}\label{prop_from_SGA6}
Let $f:A\to B$ be an injective and finite ring homomorphism
of noetherian rings. Then $f$ is the composition of a finite
sequence of injective ring homomorphisms
$$
B_n:=A\xrightarrow{f_n}B_{n-1}\xrightarrow{f_{n-1}}B_{n-2}\to
\cdots\to B_0:=B
$$
such that the induced morphism of schemes
$\Spec\,f_i:\Spec\,B_i\to B_{i+1}$ is an effective epimorphism
in the category of schemes, for every $i=1,\dots,n$.
\end{proposition}
\begin{proof} The assertion means that $f_i$ identifies
$B_{i+1}$ with the equalizer of the two natural ring homomorphisms
$(g_{i,k}:B_i\to B_i\otimes_{B_{i+1}}B_i~|~k=1,2)$, with
$g_{i,1}(b):=b\otimes 1$ and $g_{i,2}(b):=1\otimes b$ for
$i=1,\dots,n$ and every $b\in B_{i+1}$. Then let us define
inductively $B_0:=B$, and $B_{i+1}:=\Ker(h_{i,1}-h_{i,2})$ for every
$i\in\N$, where $h_{i,1},h_{i,2}:B_i\to B_i\otimes_AB_i$ are the
analogous maps.

\begin{claim}\label{cl_surge}
For every $i\in\N$, the induced map
$B_i\otimes_AB_i\to B_i\otimes_{B_{i+1}}B_i$ is an isomorphism.
\end{claim}
\begin{pfclaim} The map in question is clearly surjective.
We compute in $B_i\otimes_AB_i$ :
$$
cb\otimes b'=(c\otimes 1)\cdot(b\otimes b')=
(1\otimes c)\cdot(b\otimes b')=b\otimes cb'
\qquad
\text{for every $b,b'\in B_i$ and $c\in B_{i+1}$}
$$
which easily implies that the map in question is also injective.
\end{pfclaim}

Due to claim \ref{cl_surge}, it remains only to check that
$B_n=A$ for sufficiently large $n\in\N$. To this aim, set
$M_i:=B_i/A$, and let $S_i\subset\Spec\,A$ be the support
of $M_i$, for every $i\in\N$. Then $S_{i+1}\subset S_i$ for
every $i\in\N$, and since $A$ is noetherian, there exists
$n\in\N$ such that $S:=S_n=S_i$ for every $i\geq n$. Thus,
it suffices to check that $S=\emptyset$. If the latter
fails, let $\fp$ be a maximal point of the closed subset $S$;
a simple induction on $i\in\N$ shows that $(B_{i+1})_\fp$ is
the equalizer of the localizations
$(h_{i,1})_\fp,(h_{i,2})_\fp:B_{i,\fp}\otimes_{A_\fp}B_{i,\fp}$, for
every $i\in\N$. Hence, we may replace $f$ by its localization
$f_\fp:A_\fp\to B_\fp$, and assume from start that $(A,\fp)$ is
local, and the support of $M_n$ is the maximal ideal of $A$,
hence the $A$-module $M_n$ has finite length, so there exists
$i\geq n$ such that $M_i=M_{i+1}$, {\em i.e.} $B_i=B_{i+1}$.
Set $\kappa(\fp):=A/\fp$ and $B_i(\fp):=B_i/\fp B_i$; it
follows that each of the maps $B_i\to B_i\otimes_AB_i$ is
an isomorphism, hence the same holds for the induced maps
$B_i(\fp)\to B_i(\fp)\otimes_{\kappa(\fp)}B_i(\fp)$, so
$\dim_{\kappa(\fp)}B_i(\fp)=1$, and therefore $f:A\to B_i$
induces an isomorphism $\kappa(\fp)\isom B_i(\fp)$. By
Nakayama's lemma, we conclude that $A=B_i$, whence $M_i=0$
and $S_i=\emptyset$, a contradiction.
\end{proof}

\sset\subsubsection{}\label{subsec_perfectification}
Let $p\in\Z$ be a prime integer. For every $\F_p$-algebra $R$,
we let $R^\perf$ be the inductive limit of the system of
$\F_p$-algebras $(R_n~|~n\in\N)$ with $R_n:=R$ for every
$n\in\N$, and with transition map given by the Frobenius
endomorphism $\Phi_R:R_n\to R_{n+1}$ of $R$, for every such $n$.
The rule $R\mapsto R^\perf$ yields a left adjoint for the forgetful
functor from the category of perfect $\F_p$-algebras to the
category of all $\F_p$-algebras. Moreover, by virtue of
\cite[Th.3.5.13(ii)]{Ga-Ra}, any \'etale homomorphism $R\to S$
of $\F_p$-algebras induces an isomorphism
\set\begin{equation}\label{eq_perfect-etale}
R^\perf\otimes_RS\isom S^\perf.
\end{equation}
Especially, for every subset $\Sigma\subset R$, we have
$$
(\Sigma^{-1}R)^\perf=\Sigma^{-1}R^\perf.
$$

\begin{corollary}\label{cor_yanagihara}
Let $f:A\to B$ be a ring homomorphism.

{\em (i)}\ \
The following conditions are equivalent :
\begin{enumerate}
\alphaenu
\item
$f$ is a radicial homomorphism.
\item
$I_f\subset\cN_{B\otimes_AB}$ (notation of \eqref{subsec_notation}).
\end{enumerate}

{\em (ii)}\ \
Let moreover $p\in\N$ be a prime integer, and suppose that $A$ and
$B$ are $\F_p$-algebras, and $f$ is integral and injective. Then
{\em(a)} and {\em(b)} are also equivalent to :

\begin{enumerate}
\alphaenu\addenu\addenu
\item
$f$ induces an isomorphism $f^\perf:A^\perf\isom B^\perf$.
\end{enumerate}
\end{corollary}
\begin{proof}(a)$\Leftrightarrow$(b): Set $X:=\Spec\,A$, $Y:=\Spec\,B$,
and let $\pi:Y\times_XY\to X$ be the projection. With this notation,
$\Spec\,B\otimes_AB/I_f$ is the diagonal closed subset
$\Delta_{Y/X}\subset Y\times_XY$. Suppose then that (b) holds; it
follows already that $\Delta_{Y/X}=Y\times_XY$, so
$\phi:=\Spec\,f:Y\to X$ is injective. Moreover, let $y\in Y$ be
any point, set $x:=\phi(y)$, and let $\kappa(x)$ and $\kappa(y)$
be the residue fields of $x$ and respectively $y$; then
$Z(y):=\Spec\,\kappa(y)\otimes_{\kappa(x)}\kappa(y)=\pi^{-1}(y)$
consists of a unique point of $Y\times_XY$, which means that the
field extension $\kappa(x)\subset\kappa(y)$ is algebraic and purely
inseparable, therefore (a) holds. Conversely, if (a) holds, then
$\phi$ is injective, and $Z(y)$ consists of a unique point, so
that $\Delta_{Y/X}=Y\times_XY$, whence (b).

(b)$\Leftrightarrow$(c): Since the natural map
$R^\perf\to(R_\red)^\perf$ is an isomorphism for every $\F_p$-algebra
$R$, and since $f$ is radicial if and only if the same holds for
the induced map $f_\red:A_\red\to B_\red$, {\em we may assume that
$A$ and $B$ are reduced}, in which case (c) means that for every
$b\in B$ there exists $n\in\N$ with $b^{p^n}\in A$. Now, let us
write $B$ as the filtered colimit of the system
$(f_\lambda:A\to B_\lambda~|~\lambda\in\Lambda)$ of its finite
$A$-subalgebras. It follows easily that (c) holds for $f$ if and
only if it holds for every $f_\lambda$. On the other hand, $\Spec\,f$
is the limit of the system of maps
$(\Spec\,f_\lambda~|~\lambda\in\Lambda)$, and each $\Spec\,f_\lambda$
is surjective; it follows easily that $\Spec\,f$ is injective if
and only if the same holds for each $\Spec\,f_\lambda$. Moreover,
it is clear that for any field $K$, and any filtered system
$E_\bullet:=(E_\lambda~|~\lambda\in\Lambda)$ of field extension of $K$,
the colimit of $E_\bullet$ is algebraic and purely inseparable over
$K$ if and only if the same holds for every $E_\lambda$; summing up,
we see that $f$ is radicial if and only if the same holds for every
$f_\lambda$, hence, in view of (i), it suffices to show the equivalence
(b)$\Leftrightarrow$(c) for each $f_\lambda$, and {\em we may assume
that $f$ is a finite ring homomorphism}.

Next, let $(A_\lambda~|~\lambda\in\Lambda)$ be the filtered system of
all $\Z$-subalgebras of $A$ of finite type; also, choose any finite
subset $\Sigma\subset B$ such that $B=A[\Sigma]$. Notice that
(c) holds if and only if for every $b\in\Sigma$ there exists
$n\in\N$ such that $b^{p^n}\in A$; then the latter condition holds
if and only if there exists $\lambda\in\Lambda$ and $n\in\N$
such that $b^{p^n}\in A_\lambda$ for every $b\in\Sigma$. Moreover,
there exists a cofinal subset $\Lambda'\subset\Lambda$ such that
the restriction $f_\lambda:A_\lambda\to B_\lambda:=A_\lambda[\Sigma]$ is
a finite ring homomorphism for every $\lambda\in\Lambda'$. Furthermore,
recall that $I_f$ is the ideal generated by
$(b\otimes 1-1\otimes b~|~b\in\Sigma)$, and likewise for the ideal
$I_{f_\lambda}$, for every $\lambda\in\Lambda$. It follows that
$I_f\subset\cN_{B\otimes_AB}$ if and only if
$I_{f_\lambda}\subset\cN_{B_\lambda\otimes_{A_\lambda}B_\lambda}$ for some
$\lambda\in\Lambda'$. Summing up, it then suffices to check the
equivalence (b)$\Leftrightarrow$(c) for each map $f_\lambda$, and
therefore {\em we may assume that $A$ and $B$ are noetherian}.

Next, pick a decomposition $f=f_1\circ\cdots\circ f_n$ as in
proposition \ref{prop_from_SGA6}; since $\Spec\,f_i$ is a
surjection for every $i=1,\dots,n$, clearly $f$ is radicial
if and only if the same holds for each $f_i$. Similarly,
$f^\perf=f_1^\perf\circ\cdots\circ f_n^\perf$ is an isomorphism if
and only if $f^\perf_i$ is an isomorphism for every $i=1,\dots,n$.
Thus, it suffices to show the equivalence (b)$\Leftrightarrow$(c)
for each $f_i$, and therefore {\em we may assume that $\Spec\,f$
is an effective epimorphism of schemes}. Now, (b) says that for
every $b\in B$ there exists $n\in\N$ with
$$
(b\otimes 1-1\otimes b)^{p^n}=b^{p^n}\otimes 1-1\otimes b^{p^n}=0.
$$
By assumption, the latter holds if and only if $b^{p^n}\in A$,
whence the contention.
\end{proof}

\begin{remark}\label{rem_weak-normalization}
(i)\ \
Let $f:A\to B$ an injective ring homomorphism, and denote by
$\cE$ the set of all $A$-subalgebras $B'\subset B$ such that
the induced map $f':A\to B'$ is integral and radicial; recall
that if $B'=A[\Sigma]$ for a subset $\Sigma\subset B'$,
then $I_{f'}$ is generated by the system
$(b\otimes 1-1\otimes b~|~b\in\Sigma)$.

(ii)\ \
In view of corollary \ref{cor_yanagihara}(i), it follows
easily that if $B',B''\in\cE$, then also $B'\cdot B''\in\cE$.

(iii)\ \
Moreover, for every system $(B_\lambda~|~\lambda\in\Lambda)$
of elements of $\cE$ filtered by inclusion, we have
$\bigcup_{\lambda\in\Lambda}B_\lambda\in\cE$.

(iv)
Furthermore, if $B'\in\cE$, and $B''\subset B'$ is any
$A$-subalgebra, then $B''\in\cE$ as well.
\end{remark}

\begin{definition}
From remark \eqref{rem_weak-normalization}(ii,iii), it
follows that $\cE$ is a filtered set, and the subring
$C:=\bigcup_{B'\in\cE}B'$ is the largest element of $\cE$; we
call $C$ the {\em weak normalization} of $A$ in $B$.
\end{definition}

\begin{remark}\label{rem_funct-of-weak-norm}
(i)\ \
Consider a commutative diagram of rings :
$$
\xymatrix{ A \ar[r] \ar[d]_{g_A} & B \ar[d]^{g_B} \\
A' \ar[r] & B'
}$$
whose horizontal arrows are injective ring homomorphisms. Denote
by $C$ (resp. $C'$) the weak normalization of $A$ in $B$ (resp.
of $A'$ in $B'$). Then $g_B(C)\subset C'$. Indeed, let $x\in C$;
by remark \ref{rem_weak-normalization}(iv), the induced map
$A\to A[x]$ is radicial, hence $x\otimes 1-1\otimes x$ is
nilpotent in $A[x]\otimes_AA[x]$, so that its image
$g_B(x)\otimes 1-1\otimes g_B(x)$ is nilpotent
in $A'[x]\otimes_{A'}A'[x]$, and therefore the induced map
$A'\to A'[x]$ is radicial, whence the claim.

(ii)\ \
Let $f:A\to B$ be an injective and integral ring homomorphism;
let $j_1$ (resp. $j_2$) be the ring homomorphism
$B\to B\otimes_AB$ given by the rule : $b\mapsto b\otimes 1$
(resp. $b\mapsto 1\otimes b$) for every $b\in B$. For $i=1,2$
let also $j'_i:B\to(B\otimes_AB)_\red$ be the composition of
$j_i$ with the projection $B\otimes_AB\to(B\otimes_AB)_\red$.
Then the weak normalization of $A$ in $B$ is the equalizer of
$j'_1$ and $j'_2$. Indeed, if $b\in B$ and $A[b]$ is a radicial
$A$-algebra, then $b\otimes 1-1\otimes b$ is nilpotent in
$A[b]\otimes_AA[b]$, hence also in $B\otimes_AB$, so
$j'_1(b)=j'_2(b)$. Conversely, if the latter identity holds,
$b\otimes 1-1\otimes b$ is nilpotent in $B\otimes_AB$; but
notice that the map $A[b]\to B$ is still injective and integral,
so the kernel of the induced map $A[b]\otimes_AA[b]\to B\otimes_AB$
lies in the nilradical (lemma \ref{lem_trivialities}(iii)). It
follows that $b\otimes 1-1\otimes b$ is nilpotent in
$A[b]\otimes_AA[b]$, so $A[b]$ is a radicial $A$-algebra, as
required.
\end{remark}

\begin{lemma}\label{lem_viva-isabel}
Let $f_\bullet:=(f_\lambda:A_\lambda\to B_\lambda~|~\lambda\in\Lambda)$
be a filtered system of injective ring homomorphisms; denote
by $f:A\to B$ the colimit of $f_\bullet$, and for every
$\lambda\in\Lambda$ let $C_\lambda$ be the weak normalization
of $A_\lambda$ in $B_\lambda$. Then, the system of rings
$(B_\lambda~|~\lambda\in\Lambda)$ restricts to a filtered
system of rings $C_\bullet:=(C_\lambda~|~\lambda\in\Lambda)$,
whose colimit $C$ is the weak normalization of $A$ in $B$.
\end{lemma}
\begin{proof} Let $f':A\to C$ be the map induced by $f$, and
for every $\lambda\in\Lambda$, let
$f'_\lambda:A_\lambda\to C_\lambda$ be the map induced by $f_\lambda$.
By corollary \ref{cor_yanagihara}(i) we know that
$I_{f'_\lambda}\subset\cN_{C_\lambda\otimes_{A_\lambda}C_\lambda}$; on the
other hand, it is easily seen that $I_{f'}$ is the colimit
of the system of ideals $(I_{f'_\lambda}~|~\lambda\in\Lambda)$,
hence $I_{f'}\subset\cN_{C\otimes_AC}$. Moreover, since $C_\lambda$
is an integral $A_\lambda$-algebra for every $\lambda\in\Lambda$,
it is easily seen that $C$ is an integral $A$-algebra. Invoking
again corollary \ref{cor_yanagihara}(i), we deduce that $C$ lies
in the weak normalization $C'$ of $A$ in $B$. Conversely, if
$x\in C'$, then in particular $x$ lies in the integral closure
$C''$ of $A$ in $B$, hence we find $\lambda\in\Lambda$ and
$x_\lambda\in B_\lambda$ whose class in $B$ agrees with $x$, and
such that $x_\lambda$ is integral over $A_\lambda$; by corollary
\ref{cor_yanagihara}(i), we have
$x\otimes 1-1\otimes x\in\cN_{C'\otimes_AC'}$, and therefore
$x\otimes 1-1\otimes x\in\cN_{C''\otimes_AC''}$ as well. Hence, after
replacing $\lambda$ by some $\lambda'\geq\lambda$, we may assume
that $x_\lambda\otimes 1-1\otimes x_\lambda\in
\cN_{C''_\lambda\otimes_{A_\lambda}C''_\lambda}$, where $C''_\lambda$ denotes
the integral closure of $A_\lambda$ in $B_\lambda$. Now, the inclusion
$A_\lambda[x]\to C''_\lambda$ is integral and injective, hence
the kernel of the induced ring homomorphism
$A_\lambda[x]\otimes_{A_\lambda}A_\lambda[x]\to
C''_\lambda\otimes_{A_\lambda}C''_\lambda$ lies in the nilradical
(lemma \ref{lem_trivialities}(iii)). Summing up, we conclude that
$x_\lambda\otimes 1-1\otimes x_\lambda$ is already nilpotent in
$A_\lambda[x]\otimes_{A_\lambda}A_\lambda[x]$, so $x_\lambda\in C_\lambda$,
and finally $x\in C$. This shows that $C'=C$, as required.
\end{proof}

\begin{proposition}\label{prop_weak-norm-and-smooth-bc}
Let $f:A\to B$ be an injective ring homomorphism, $C$ the weak
normalization of $A$ in $B$, and $A'$ a flat $A$-algebra such
that the induced morphism $\Spec\,A'\to\Spec\,A$ has geometrically
reduced fibres. Suppose that :
\begin{itemize}
\item
either $A'$ is a finitely presented $A$-algebra
\item
or else, $A$ is noetherian.
\end{itemize}
Then $A'\otimes_Af:A'\to B':=A'\otimes_AB$ is an injective ring
homomorphism, and the weak normalization of $A'$ in $B'$ is
$A'\otimes_AC$.
\end{proposition}
\begin{proof}Since $A'$ is a flat $A$-algebra, $A'\otimes_Af$ is
still injective, as stated. Now, let $A^\nu$ (resp. $A'^\nu$) be
the integral closure of $A$ in $B$ (resp. of $A'$ in $B'$); then
$A'^\nu=A'\otimes_AA^\nu$ (proposition \ref{prop_normal-and-geom-red-bc}).
We may then replace $B$ by $A^\nu$, and assume that $B$ is an
integral $A$-algebra. In this case, notice that for every
$A$-algebra $R$ we have $(A'\otimes_AR)_\red=A'\otimes_AR_\red$
(proposition \ref{prop_normal-and-geom-red-bc}); then the assertion
follows easily from remark \ref{rem_funct-of-weak-norm}(ii).
\end{proof}

\sset\subsubsection{}\label{subsec_pflieger}
Let $f:A\to B$ be an integral and injective ring homomorphism,
$I\subset A$ an ideal; set $X:=\Spec\,A$, $Y:=\Spec\,B$,
$Z:=\Spec\,A/I$, $U:=X\setminus Z$, and notice that
$\phi:=\Spec\,f:Y\to X$ is surjective by the going up theorem,
hence the same holds for its restriction $\phi^{-1}Z\to Z$, and
therefore the kernel of $A/I\otimes_Af:A/I\to B/IB$ lies in the
nilradical (lemma \ref{lem_trivialities}(ii)), {\em i.e.} the
induced map $(A/I)_\red\to(B/IB)_\red$ is injective and integral
as well.

\begin{lemma}\label{lem_pflieger}
In the situation of \eqref{subsec_pflieger}, suppose that $\phi$
restricts to a radicial morphism of schemes
$\phi_{|U}:\phi^{-1}U\to U$, and let $C$ (resp. $\bar C$) be the
weak normalization of $A$ in $B$ (resp. of $(A/I)_\red$ in
$(B/IB)_\red$). Then we have :
$$
C=B\times_{(B/IB)_\red}\bar C.
$$
\end{lemma}
\begin{proof} By remark \ref{rem_funct-of-weak-norm}(i), the
projection $B\to(B/IB)_\red$ maps $C$ into $\bar C$, whence
the inclusion $C\subset C':=B\times_{(B/IB)_\red}\bar C$. For
the converse inclusion, we need to show that the induced
map $g:A\to C'$ is radicial. To this aim, it suffices to
check that the same holds for $A/I\otimes_Ag:A/I\to C'/IC$
and for $A[a^{-1}]\otimes_Ag:A[a^{-1}]\to C'[a^{-1}]$, for every
$a\in I$. However, for any such $a$, the localization $C'[a^{-1}]$
is the fibre product $B[a^{-1}]\times_{(B/IB)_\red[a^{-1}]}\bar C[a^{-1}]$,
and since $(B/IB)_\red[a^{-1}]=\bar C[a^{-1}]=0$, we see that
$C'[a^{-1}]=B[a^{-1}]$; but $f$ induces a radicial map
$A[a^{-1}]\to B[a^{-1}]$, since $\phi_{|U}$ is radicial. It
remains thus only to check that $C/IC$ is a radicial $A/I$-algebra;
to this aim, notice that the kernel of the surjection
$B\to(B/IB)_\red$ is the radical $J$ of the ideal $IB\subset B$.
Then $J$ is also the kernel of the surjection $C'\to\bar C$;
on the other hand, since $IB\subset J\subset C$, we have
$I^2B\subset IC\subset IB\subset J$, and clearly $J$ is also
the radical of the ideal $I^2B$ of $C$. Hence, $C/IC$ is a
radicial $A/I$-algebra if and only if the same holds for the
$A/I$-algebra $C/J$; but $C/J=\bar C$, whence the contention.
\end{proof}

\begin{definition}\label{def_p-integrally-closed}
(i)\ \
Let $R$ be a ring, $S\subset R$ a subring, and $p>0$
a prime integer. We say that {\em $S$ is $p$-integrally
closed in $R$}, if for every $x\in R$ such that $x^p\in S$,
we have $x\in S$.

(ii)\ \
Let $R$ be a ring and $p>0$ a prime integer; if
$(S_\lambda~|~\lambda\in\Lambda)$ is any family of
$p$-integrally closed subrings of $R$, then clearly
$\bigcap_{\lambda\in\Lambda}S_\lambda$ is also $p$-integrally
closed in $R$. Hence, for every subring $S\subset R$
there exists a smallest $p$-integrally closed subring
$S'$ of $R$ containing $S$. We say that $S'$ is the
{\em $p$-integral closure of $S$ in $R$}.
\end{definition}

\begin{remark}\label{rem_charact-p-int-clos}
(i)\ \
Let $R$ be a ring, $S\subset R$ a subring; consider the
increasing countable system $(S_n~|~n\in\N)$ of subrings
of $R$ defined inductively as follows : $S_0:=S$, and
$S_{n+1}:=S_n[\Sigma_n]$ for every $n\in\N$, where
$\Sigma_n:=\{x\in R~|~x^p\in S_n\}$. Then it is easily
seen that $\bigcup_{n\in\N}S_n$ is the $p$-integral
closure of $S$ in $R$.

(ii)\ \
In the situation of remark \ref{rem_funct-of-weak-norm}(i),
let $D$ (resp. $D'$) be the $p$-integral closure of $A$ in
$B$ (resp. of $A'$ in $B'$). Define the system of subrings
$(A_n~|~n\in\N)$ of $B$ (resp. $(A'_n~|~n\in\N)$ of $B'$)
as in (i); a simple induction shows that $g_B(A_n)\subset A'_n$
for every $n\in\N$, so $g_B$ restricts to a ring homomorphism
$D\to D'$.
\end{remark}

\begin{lemma}\label{lem_crit-p-integr-closed}
{\em(i)}\ \
Let $R$ be a ring, $a\in R$ a regular element such that
$pR\subset a^pR$. Then $R$ is $p$-integrally closed in
$R[a^{-1}]$ if and only if the Frobenius endomorphism of
$R/a^pR$ induces an injective ring homomorphism
$\bar\Phi:R/aR\to R/a^pR$. 

{\em(ii)}\ \
Let $R$ be a topological ring, $S\subset R$ an open subring,
$T$ the $p$-integral closure of $S$ in $R$; endow $S$ and $T$
with the topologies induced by $R$, and denote by $R^\wedge$,
$S^\wedge$ and $T^\wedge$ the respective completions. Then
$T^\wedge$ is the $p$-integral closure in $R^\wedge$ of $S^\wedge$.

{\em(iii)}\ \
Let $R_\bullet:=(R_\lambda~|~\lambda\in\Lambda)$ be a system of
rings, with transition morphisms $f_{\lambda\mu}:R_\lambda\to R_\mu$
for every $\lambda,\mu\in\Lambda$ with $\mu\geq\lambda$. For
every $\lambda\in\Lambda$, let also $S_\lambda\subset R_\lambda$
be a subring, such that $f_{\lambda\mu}(S_\lambda)\subset S_\mu$ for
every such $\lambda,\mu$. Denote by $D_\lambda\subset R_\lambda$
the $p$-integral closure of $S_\lambda$ in $R_\lambda$, for every
$\lambda\in\Lambda$. The following holds :
\begin{enumerate}
\alphaenu
\item
Let $R$ and $S$ be the limits of $R_\bullet$ and respectively
of $S_\bullet:=(S_\lambda~|~\lambda\in\Lambda)$. Then the limit
of the induced system $D_\bullet:=(D_\lambda~|~\lambda\in\Lambda)$
is the $p$-integral closure of $S$ in $R$.
\item
Suppose moreover that $\Lambda$ is filtered, and let $R'$ and
$S'$ be the colimits of $R_\bullet$ and respectively of $S_\bullet$.
Then the colimit of $D_\bullet$ is the $p$-integral closure of
$S'$ in $R'$.
\end{enumerate}

{\em(iv)}\ \
Let $B$ be a ring, $A\subset B$ a subring, $C$ the $p$è-integral
closure of $A$ in $B$, and $S\subset A$ a multiplicative subset.
Then $S^{-1}C$ is the $p$-integral closure of $S^{-1}A$ in $S^{-1}B$.
\end{lemma}
\begin{proof}(i): Suppose first that $\bar\Phi$ is injective;
let $x\in R[a^{-1}]$ be any element such that $x^p\in R$, and
suppose, by way of contradiction, that $x\notin R$. There exists
a smallest $m\in\N$ such that $y:=a^mx\notin R$ and $a^{m+1}x\in R$.
Hence, $y^p=a^{pm}x^p\in R$, and therefore
$(a^{m+1}x)^p=a^py^p\in a^pR$. Thus, the class of $a^{m+1}x$
lies in $\Ker\,\bar\Phi$, {\em i.e.} $a^{m+1}x=az$ for some
$z\in R$; since $a$ is regular, it follows that $z=a^mx\in R$,
contradicting the choice of $m$.

Conversely, suppose that $R$ is integrally closed in $R[a^{-1}]$,
and let $x\in R$ be an element whose class $\bar x\in R/aR$ lies
in $\Ker\,\bar\Phi$; this means that there exists $y\in R$ such
that $x^p=a^py$. Hence, $(x/a)^p\in R$, so that $x/a\in R$,
and therefore $\bar x=0$.

(ii): Notice that $S^\wedge\cdot T=T^\wedge$, hence $T^\wedge$
lies in the $p$-integral closure of $S^\wedge$ in $R^\wedge$.
Hence, we are reduced to checking that $S$ is $p$-integrally
closed in $R$ if and only if $S^\wedge$ is $p$-integrally closed
in $R^\wedge$.
Suppose first that $S$ is $p$-integrally closed in $R$,
and let $x\in R^\wedge$ be any element such that $x^p\in S^\wedge$.
Denote by $j:R\to R^\wedge$ the completion map.
Pick a Cauchy net $(x_\lambda~|~\lambda\in\Lambda)$ in $R$
whose limit is $x$; then $(x_\lambda^p~|~\lambda\in\Lambda)$
is a Cauchy net converging to $x^p$, and since $S^\wedge$ is
open in $R^\wedge$ (corollary \ref{cor_not-in-Bourbaki}(i)),
we may replace $\Lambda$ by a cofinal subset, and assume that
$j(x_\lambda)^p\in S^\wedge$ for every $\lambda\in\Lambda$. Then
$x_\lambda^p\in S$ (corollary \ref{cor_not-in-Bourbaki}(ii)),
whence $x_\lambda\in S$ for every $\lambda\in\Lambda$, and
finally $x\in S^\wedge$, as required. Conversely, if $S^\wedge$
is $p$-integrally closed in $R^\wedge$, and $x\in R$ is any
element with $x^p\in S$, we have $j(x)^p\in S^\wedge$, so
$j(x)\in S$, whence $x\in S$, again by corollary
\ref{cor_not-in-Bourbaki}(ii).

The proof of (iii) shall be left to the reader.

(iv): Let $b\in B,s\in S$ such that $(b/s)^p\in S^{-1}C$; it follows
easily that $t^pb^p\in C$ for some $t\in S$, hence $tb\in C$, and
consequently $b/s\in S^{-1}C$. This shows that $S^{-1}C$ is
$p$-integrally closed in $S^{-1}B$, hence it contains the
$p$-integral closure $C'$ of $S^{-1}A$ in $S^{-1}B$. On the
other hand, by remark \ref{rem_charact-p-int-clos}, the localization
$B\to S^{-1}B$ restricts to a ring homomorphism $C\to C'$; since
$C'$ is an $S^{-1}A$-algebra, we must then have $S^{-1}C\subset C'$.
\end{proof}

\sset\subsubsection{}\label{subsec_global-weak-norm}
Let $X$ be a scheme, $\cA\to\cB$ a monomorphism of quasi-coherent
$\cO_X$-algebras; for every affine open subset $U\subset X$, let
$\cC(U)$ (resp. $\cD(U)$) be the weak normalization (resp. the
$p$-integral closure) of $\cA(U)$ in $\cB(U)$; by remark
\ref{rem_funct-of-weak-norm}(i), the rule : $U\mapsto\cC(U)$
(resp. $U\mapsto\cD(U)$) defines a subpresheaf of $\cO_X$-algebras
$\cC\subset\cB$ (resp. $\cD\subset\cB$) on the site of affine open
subsets of $X$, and it follows easily from proposition
\ref{prop_weak-norm-and-smooth-bc} (resp. from lemma
\ref{lem_crit-p-integr-closed}(iv)) that $\cC$ (resp. $\cD$)
is a sheaf on this site, so it extends uniquely to a subsheaf
of $\cB$ on the Zariski site of $X$; moreover, $\cC$ and $\cD$
are quasi-coherent $\cO_X$-algebras. We call $\cC$ and $\cD$
respectively the {\em weak normalization} and the
{\em $p$-integral closure} of $\cA$ in $\cB$.

\begin{proposition}\label{prop_many-rings-equal}
Let $f:A\to B$ be an injective ring homomorphism, $p\in\N$ a
prime integer, and suppose that $f$ induces an isomorphism
$A[1/p]\isom B[1/p]$. Let $A^\nu$ (resp. $C_0$, resp. $C_1$) be
the integral closure (resp. the weak normalization, resp. the
$p$-integral closure) of $A$ in $B$; denote by
$(A/pA)_\red\xrightarrow{j}(A^\nu/pA^\nu)_\red\xleftarrow{\pi}A^\nu$
the natural maps, and set
$$
\begin{aligned}
C_2:=&\{b\in B~|~\text{there exists $n\in\N$ such that $b^{p^n}\in A$}\} \\
D:=&\{x\in(A^\nu/pA^\nu)_\red~|~
\text{there exists $n\in\N$ such that $x^{p^n}\in\Img(j)$}\}.
\end{aligned}
$$
Then $j$ is injective, and $C_0=C_1=C_2=\pi^{-1}D$.
\end{proposition}
\begin{proof} One argues as in \eqref{subsec_pflieger} to see that
$j$ is injective.

Now, the inclusion $C_2\subset C_1$ is obvious, and the identity
$C_0=\pi^{-1}D$ follows easily from lemma \ref{lem_pflieger}. By
remark \ref{rem_charact-p-int-clos}, we may write $C_1$ as the
increasing union of the family
$A_0\subset A_1\subset A_2\subset\cdots$ of subrings of $B$,
defined inductively by $A_0:=A$, and $A_{i+1}:=A_i[\Sigma_i]$
for every $i\in\N$, where $\Sigma_i:=\{b\in B~|~b^p\in A_i\}$.
By remark \ref{rem_weak-normalization}(iii), in order to check
that $C_1\subset C_0$, it will therefore suffice to show that the
inclusion map $A_i\to A_{i+1}$ is radicial for every $i\in\N$.
However, notice that $A_i[1/p]=A_{i+1}[1/p]$ for every such $i$,
hence we are reduced to checking that the induced map
$A_i/pA_i\to A_{i+1}/pA_{i+1}$ is radicial, which is clear by
construction.

Lastly, let us show that $C_0\subset C_2$. To this aim, let us write
$f$ as the colimit of a filtered system of injective ring homomorphisms
$(f_\lambda:A_\lambda\to B_\lambda~|~\lambda\in\Lambda)$, where
$A_\lambda$ and $B_\lambda$ are $\Z$-subalgebras of finite type of $A$
and respectively $B$, for every $\lambda\in\Lambda$; since $f$
induces an isomorphism $A[1/p]\isom B[1/p]$, we may also assume that
$f_\lambda$ induces an isomorphism $A_\lambda[1/p]\isom B_\lambda[1/p]$
for every $\lambda\in\Lambda$ (details left to the reader). For every
such $\lambda$, let $C_{2,\lambda}$ be the subset of all $b\in B_\lambda$
such that $b^{p^n}\in A_\lambda$ for some $n\in\N$, and let $C_{0,\lambda}$
be the weak normalization of $A_\lambda$ in $B_\lambda$. Clearly $C_2$
is the filtered union of the resulting system
$(C_{2,\lambda}~|~\lambda\in\Lambda)$; on the other hand, $C_0$ is the
colimit of the system of subrings $(C_{0,\lambda}~|~\lambda\in\Lambda)$,
by lemma \ref{lem_viva-isabel}. Hence, we are reduced to checking
that $C_{0,\lambda}\subset C_{2,\lambda}$ for every $\lambda\in\Lambda$,
and {\em we may assume that $A$ and $B$ are noetherian}.

Let now $b\in C_0$, and set $B':=A[b]$; we have to check that
$b^{p^n}\in A$ for some $n\in\N$, and obviously $A[1/p]=B'[1/p]$,
so we may replace $B$ by $B'$ and {\em suppose that $f$ is finite
and radicial}. In this case, pick a decomposition
$f=f_1\circ\cdots\circ f_n$ as in proposition \ref{prop_from_SGA6};
again, it suffices to check the assertion for each $f_i$, hence
{\em we may assume that $\Spec\,f$ is an effective epimorphism of
schemes} (and still radicial). Thus, for every $b\in B$ the element
$c:=b\otimes 1-1\otimes b$ is nilpotent in $B\otimes_AB$ (corollary
\ref{cor_yanagihara}(i)), and moreover notice that the two natural
maps $B[1/p]\to B\otimes_AB[1/p]$ are isomorphisms, since
$A[1/p]=B[1/p]$, hence the same holds for the multiplication law
$B\otimes_AB[1/p]\to B[1/p]$, and the therefore image of $c$ vanishes
in $B\otimes_AB[1/p]$, {\em i.e.} $p^kc=0$ in $B\otimes_AB$ for
some $k\in\N$, and after replacing $k$ by a larger integer, we
may assume that $c^k=0$ as well. We write
$$
b^{p^n}\otimes 1=1\otimes b^{p^n}+
\sum_{j=1}^{p^n}\binom{p^n}{j}(1\otimes b^{p^n-j})\cdot c^j=
1\otimes b^{p^n}+\sum_{j=1}^{k-1}\binom{p^n}{j}(1\otimes b^{p^n-j})\cdot c^j
$$
and recall that the $p$-adic valuation of $\binom{p^n}{j}$ equals
$n-v_p(j)$ for $j=1,\dots,p^n$, where $v_p(j)$ denotes the $p$-adic
valuation of $j$. In particular, with $n:=2k$, we see that
$\binom{p^n}{j}$ is divisible by $p^k$ for every $j=1\dots,k-1$,
so finally $b^{p^n}\otimes 1=1\otimes b^{p^n}$ in $B\otimes_AB$,
and consequently $b^{p^n}\in A$, since $\Spec\,f$ is an effective
epimorphism.
\end{proof}

\begin{proposition}\label{prop_invariance-under-almiso}
Let us resume the situation of remark
{\em\ref{rem_funct-of-weak-norm}(i)}, and let $a\in A$ be any element
and $n\in\N$ any integer; also, denote by $A^\nu$ (resp. $A'^\nu$)
the integral closure of $A$ in $B$ (resp. of $A'$ in $B'$). The
following holds :

{\em (i)}\ \
If $g_A$ and $g_B$ are injective, $aB'\subset\Img\,g_B$
and $a^nA'\subset\Img\,g_A$, then :
$$
aA'^\nu\subset g_B(A^\nu)
\qquad\text{and}\qquad
aC'\subset g_B(C).
$$

{\em (ii)}\ \
If $g_A$ and $g_B$ are surjective, and $a\Ker\,g_B=0$, then :
$$
a\cdot g_B^{-1}A'^\nu\subset A^\nu
\qquad\text{and}\qquad
a\cdot g_B^{-1}C'\subset C.
$$

{\em (iii)}\ \
If $g_A$ and $g_B$ are surjective, and $\Ker\,g_B\subset\cN_B$
(notation of \eqref{subsec_notation}), then :
$$
g_B^{-1}A'^\nu=A^\nu
\qquad\text{and}\qquad
g_B^{-1}C'=C.
$$
\end{proposition}
\begin{proof}(i): Let $b\in A'^\nu$, and pick any monic polynomial
$P:=X^m+a'_1X^{m-1}+\cdots+a'_m\in A'[X]$ such that $P(b^n)=0$. Hence
$a^{mn}\cdot P(b^n)=0=(ab)^{mn}+a^na'_1(ab)^{n(m-1)}+\cdots+a^{mn}a'_m=0$,
and by assumption $a^na'_1,a^{2n}a'_2,\dots,a^{mn}a'_m\in g_A(A)$ and
$ab\in g_B(B)$, so $ab\in g_B(A^\nu)$, whence the first stated
inclusion. Next, set $C_0:=g_A(A)+aC'$; clearly $C'_0$ is an
$A$-subalgebra of $g_B(B)$, and we need to check that the induced
map $A\to C_0$ is radicial.

\begin{claim}\label{cl_viva-isabel}
The kernel and cokernel of the induced map
$C_0\otimes_AC_0\to C'\otimes_{A'}C'$ are annihilated by $a^{n+2}$.
\end{claim}
\begin{pfclaim} By consruction, the cokernel of the inclusion
$j:C_0\to C'$ is annihilated by $a$; since the kernel of
$C'\otimes_Aj$ is a quotient of $\Tor_1^A(C',C'/C_0)$, it follows
that both the kernel and cokernel of $C'\otimes_Aj$ are annihilated
by $a$. By the same token, the kernel and cokernel of
$(C'\otimes_AC_0)\otimes_{C_0}j:C'\otimes_AC_0\to C'\otimes_AC'$ are
annihilated by $a$ as well. Lastly, the natural map
$\pi:C'\otimes_AC'\to C''\otimes_{A'}C'$ is clearly surjective, and
its kernel is generated by the elements of the form
$a'x\otimes y-x\otimes a'y$, for all $a'\in A'$ and all $x,y\in C'$.
By assumption, $a^n(a'x\otimes y-x\otimes a'y)=0$ in $C'\otimes_AC'$,
so $a^n\Ker\,\pi=0$. The claim follows easily.
\end{pfclaim}

Let $c\in aC'$ be any element; then $c\otimes 1-1\otimes c$ is
nilpotent in $C'\otimes_{A'}C'$ (corollary \ref{cor_yanagihara}(i));
by claim \ref{cl_viva-isabel}, it follows that there exists $N\in\N$
with $a^{n+2}(c\otimes 1-1\otimes c)^N=0$ in $C_0\otimes_AC_0$. Then,
say that $c=ad$ for some $d\in C'$; it follows that
$(c\otimes 1-1\otimes c)^3=a^3d^3\otimes 1-3a^2d^2\otimes ad+
3ad\otimes a^2d^2+1\otimes a^3d^3$ is divisible by $a$ in
$C_0\otimes_AC_0$, and finally $(c\otimes 1-1\otimes c)^{N+3n+6}=0$
in $C_0\otimes_AC_0$, whence the contention.

(ii): Let $b$ and $P$ be as in the foregoing; suppose that
$b=g_B(c)$ for some $c\in B$, and pick $a_1,\dots,a_m\in A$
with $g_A(a_i)=a'_i$ for $i=1,\dots,m$. Set
$Q:=X^m+a_1X^{m-1}+\cdots+a_m\in A[X]$; then $aQ(c)=0$, and it
follows easily that $ac\in A^\nu$, whence the first inclusion
of (ii). Next, set $C_0:=A+a\cdot g_B^{-1}C'$; we need to show
that the resulting ring homomorphism $A\to C_0$ is radicial.
However, notice that the restriction $g_C:C_0\to C$ of $g_B$
has both kernel and cokernel annihilated by $a$; then by lemma
\ref{lem_extract-afg-afp} we find an $A$-linear map $h:C\to C_0$
with $h\circ g_C=a^2\one_{C_0}$ and $g_C\circ h=a^2\cdot\one_C$.
Hence $(h\otimes_Ah)\circ(g_C\otimes_Ag_C)=a^4\one_{C_0}$ and
$(g_C\otimes_Ag_C)\circ(h\otimes_Ah)=a^4\one_C$, so both the
kernel and cokernel of
$g_C\otimes_Ag_C:C_0\otimes_AC_0\to C\otimes_AC=C\otimes_{A'}C$
are annihilated by $a^4$. We may now argue as in the proof of
(i), to conclude that $c\otimes 1-1\otimes c$ is nilpotent in
$C_0\otimes_AC_0$, for every $c\in a\cdot g_B^{-1}C'$, whence the
contention.

(iii): Let $b,c,P$ and $Q$ be as in the proof of (ii); since
$\Ker\,g_B\subset\cN_B$, we have $Q(c)^n=0$ for some $n\in\N$,
so $c\in A^\nu$, whence the first identity of (iii). Next, set
$C_0:=g_B^{-1}C'$; we need to check that the map $A\to C_0$ is
radicial, or equivalently, that the same holds the induced map
on reduced quotients $A_\red\to(C_0)_\red$. However, $g_B$
induces an isomorphism $(C_0)_\red\isom(C')_\red$, and by
definition the map $A'_\red\to(C')_\red$ is radicial; since
$g_A$ induces a surjection $A_\red\to A'_\red$, the assertion
follows easily.
\end{proof}

\section{Cohomology and local cohomology of sheaves}
\label{chap_local-coh}

\subsection{Cohomology of topoi and topological spaces}
\begin{definition}
Let $X$ be a topological space, $\cF$ a sheaf on $X$.
\begin{enumerate}
\item
We say that $\cF$ is {\em flabby\/} if the restriction map
$\cF(X)\to\cF(U)$ is surjective for every open subset $U\subset X$.
\item
We say that $\cF$ is {\em qc-flabby\/} if the restriction
map $\cF(V)\to\cF(U)$ is a surjection whenever $U\subset V$
are quasi-compact open subsets of $X$.
\item
The {\em support\/} of $\cF$ is the subset :
$$
\Supp\,\cF:=\{x\in X~|~\cF_x\neq\emptyset\}\subset X.
$$
\item
On the other hand, if $\cF$ is an abelian sheaf, we
define the {\em support} of $\cF$ as the subset :
$$
\Supp\,\cF:=\{x\in X~|~\cF_x\neq 0\}\subset X.
$$
\end{enumerate}
\end{definition}

\begin{lemma}\label{lem_flabby-is-local}
Let $X$ be a topological space, $(X_\lambda~|~\lambda\in\Lambda)$ a
family of open subsets of $X$ with $X=\bigcup_{\lambda\in\Lambda}X_\lambda$,
and $\cF$ a sheaf on $X$. For every $\lambda\in\Lambda$, let
$\cF_\lambda$ be the restriction of $\cF$ to $X_\lambda$. Then $\cF$
is flabby if and only if $\cF_\lambda$ is flabby for every
$\lambda\in\Lambda$.
\end{lemma}
\begin{proof} It is easily seen that if $\cF$ is flabby, the
same holds for every $\cF_\lambda$. Conversely, suppose that
$\cF_\lambda$ is flabby for every $\lambda\in\Lambda$; let
$U\subset X$ be any open subset, and $\sigma\in\cF(U)$. For
every subset $\Lambda'\subset\Lambda$ set
$X_{\Lambda'}:=\bigcup_{\lambda\in\Lambda'}X_\lambda$, and denote
by $\cP$ the set of all pairs $(\Lambda',\sigma')$ where
$\Lambda'$ is a subset of $\Lambda$, and
$\sigma'\in\cF(U\cup X_{\Lambda'})$ is any section whose
restriction to $U$ agrees with $\sigma$. We endow $\cP$ with
the partial order such that
$(\Lambda',\sigma')\leq(\Lambda'',\sigma'')$ if and only if
$\Lambda'\subset\Lambda''$ and the restriction of $\sigma''$
to $U\cup X_{\Lambda'}$ agrees with $\sigma'$. It is clear
that if $((\Lambda_i,\sigma_i)~|~i\in I)$ is any totally
ordered subset in $\cP$, then there exists $(\Lambda',\sigma')$
in $\cP$ with $\Lambda'=\bigcup_{i\in I}\Lambda_i$ and such that
$(\Lambda_i,\sigma_i)\leq(\Lambda',\sigma')$ for every $i\in I$.
Also, $(\emptyset,\sigma)\in\cP$. By Zorn's lemma, $\cP$ admits
then a maximal element $(\Lambda',\sigma')$, and it suffices to
check that $\Lambda'=\Lambda$. Thus, suppose that
$\lambda\in\Lambda\setminus\Lambda'$; set
$\Lambda'':=\Lambda'\cup\{\lambda\}$, and let $\sigma''$ be the
restriction of $\sigma'$ to $U_\lambda:=X_\lambda\cap(U\cup X_{\Lambda'})$.
Since $\cF_\lambda$ is flabby, there exists $\sigma_\lambda\in X_\lambda$
whose restriction to $U_\lambda$ agrees with $\sigma''$. Then there
exists $\tau\in\cF(U\cup X_{\Lambda''})$ whose restriction
to $U\cup X_{\Lambda'}$ agrees with $\sigma'$ and whose restriction
to $X_\lambda$ agrees with $\sigma_\lambda$. By construction, we have
$(\Lambda'',\tau)>(\Lambda',\sigma')$, a contradiction.
\end{proof}

\sset\subsubsection{}\label{subsec_flabby-short-exact}
Let $X$ be a topological space, and consider a short
exact sequence
$$
\Sigma\qquad :\qquad
0\to\cF'\to\cF\to\cF''\to 0
\qquad\qquad\qquad
$$
of abelian sheaves on $X$. For every open subset $U\subset X$,
the sequence $\Sigma$ induces a complex :
$$
\Sigma(U)\qquad :\qquad
0\to\cF'(U)\to\cF(U)\to\cF''(U)\to 0.
\qquad\qquad\qquad
$$

\begin{lemma}\label{lem_qc-flabby}
In the situation of \eqref{subsec_flabby-short-exact}, the
following holds :
\begin{enumerate}
\item
If $\cF'$ is flabby, the sequence $\Sigma(U)$ is short exact.
\item
If both $\cF$ and $\cF'$ are flabby, the same holds for $\cF''$.
\item
Suppose that $X$ is locally coherent and quasi-separated. Then :
\begin{enumerate}
\item
If $\cF'$ is qc-flabby and $U$ is quasi-compact, the sequence
$\Sigma(U)$ is short exact.
\item
If both $\cF$ and $\cF'$ are qc-flabby, the same holds for $\cF''$.
\end{enumerate}
\end{enumerate}
\end{lemma}
\begin{proof}(i): Let $s''\in\cF''(U)$ be any section; we need
to show that $s''$ is the image of an element of $\cF(U)$. To
this aim, let $\cP$ be the set of all pairs $(V,t)$, where
$V\subset U$ is an open subset, and $t\in\cF(V)$ is a section
whose image in $\cF''(V)$ agrees with the restriction of $s''$.
Then $\cP$ is partially ordered by declaring that $(V,t)\leq(V',t')$
if and only if $V\subset V'$ and $t$ is the restriction of $t'$
to $V$. Clearly, if $((V_\lambda,t_\lambda)~|~\lambda\in\Lambda)$
is a totally ordered subset of $\cP$, then there exists a unique
section $t$ on $V:=\bigcup_{\lambda\in\Lambda}V_\lambda$ such that
$(V,t)\in\cP$ and $(V_\lambda,t_\lambda)\leq(V,t)$ for every
$\lambda\in\Lambda$. Moreover, since $\cF\to\cF''$ is an
epimorphism, it is easily seen that $\cP\neq\emptyset$. By
Zorn's lemma, $\cP$ admits therefore a maximal element $(V,s)$,
and it suffices to check that $V=U$. Thus, suppose that
$x\in U\setminus V$; by assumption, we may find an open
neighborhood $V'$ of $x$ in $U$, and $t\in\cF(V')$
with $(V',t)\in\cP$. Set $W:=V\cap V'$, and let $t_W$
(resp. $s_W$) be the restriction of $t$ (resp. of $s$)
to $W$; then $s_W-t_W\in\cF'(W)$. Since $\cF'$ is flabby,
the same holds for its restriction $\cF'_{|V}$ (lemma
\ref{lem_flabby-is-local}), hence we find $r\in\cF'(V')$
whose restriction to $W$ agrees with $s_W-t_W$; set
$t':=r+t\in\cF(V')$. Set $V'':=V\cup V'$; clearly the
restriction of $t'$ to $W$ agrees with $s_W$, hence there
exists $u\in\cF(V'')$ whose restriction to $V$ and $V'$
agrees respectively with $s$ and $t'$. By construction,
we have $(V'',t')\in\cP$ and $(V'',t')>(V,s)$, a contradiction.

(ii): In view of (i), we get a commutative ladder with
short exact rows :
$$
\xymatrix{ 0 \ar[r] & \cF'(X) \ar[r] \ar[d]_{\rho'} &
\cF(X) \ar[r] \ar[d]^\rho & \cF''(X) \ar[r] \ar[d]^{\rho''} & 0 \\
0 \ar[r] & \cF'(U) \ar[r] & \cF(U) \ar[r] &
\cF''(U) \ar[r] & 0
}$$
whence, by the snake lemma, a surjection :
$\Coker\,\rho\to\Coker\,\rho''$. But $\Coker\,\rho=0$, since
$\cF$ is flabby; so $\Coker\,\rho''=0$, whence the assertion.

The proof of (iii.a) is a variant of that of (i). Indeed,
we may replace $X$ by $U$, and thereby assume that $X$ is
quasi-compact. Then we have to check that every section
$s''\in\cF''(X)$ is the image of an element of $\cF(X)$.
However, we can find a finite covering $(U_i~|~i=1,\dots,n)$
of $X$, consisting of quasi-compact open subsets, such that
$s''_{|U_i}$ is the image of a section $s_i\in\cF(U_i)$. For
every $k\leq n$, let $V_k:=U_1\cup\cdots\cup U_k$;
we show by induction on $k$ that $s''_{|V_k}$ is in the
image of $\cF(V_k)$; the lemma will follow for $k=n$.
For $k=1$ there is nothing to prove. Suppose that $k>1$
and that the assertion is known for all $j<k$; hence
we can find $t\in\cF(V_{k-1})$ whose image is $s''_{|V_{k-1}}$.
The difference $u:=(t-s_k)_{|U_k\cap V_{k-1}}$ lies in the
image of $\cF'(U_k\cap V_{k-1})$; since $U_k\cap V_{k-1}$
is quasi-compact and $\cF'$ is qc-flabby, $u$ extends to
a section of $\cF'(U_k)$. We can then replace $s_k$ by
$s_k+u$, and assume that $s_k$ and $t$ agree on $U_k\cap V_{k-1}$,
whence a section on $V_k=U_k\cup V_{k-1}$ with the sought
property. Assertion (iii.b) is proved as (ii).
\end{proof}

\begin{lemma}\label{lem_franziska}
Let $f:Y\to X$ be a continuous map of topological spaces,
$\cF$ a sheaf on $X$, and $\cG$ a sheaf on $Y$. We have :
\begin{enumerate}
\item
If $f$ is an open immersion and $\cF$ is qc-flabby
on $X$, then $f^*\cF$ is qc-flabby on $Y$.
\item
If\/ $\cG$ is flabby on $Y$, then $f_*\cG$ is flabby on $X$,
and $R^pf_*\cG=0$ for every $p>0$.
\item
If $f$ is quasi-compact and $\cG$ is qc-flabby on $Y$, then
$f_*\cG$ is qc-flabby on $X$.
\item
If $X$ and $Y$ are locally coherent, $f$ is
quasi-compact and quasi-separated, and $\cG$ is
abelian and qc-flabby, then $R^pf_*\cG=0$ for every $p>0$.
\item
Let $\cA$ be any sheaf of rings on $X$. Every injective
$\cA$-module is flabby.
\item
If $X$ is locally noetherian, every qc-flabby sheaf
on $X$ is flabby.
\end{enumerate}
\end{lemma}
\begin{proof}(i) and (iii) are trivial.

(ii): Clearly $f_*\cG$ is flabby, and it remains to check
that $R^pf_*\cG=0$ for every $p>0$. To this aim, in view
of remark \ref{rem_acyclic-crit}(iv,v,vi) and lemma
\ref{lem_qc-flabby}(ii), it suffices to show that for
every short exact sequence $0\to\cG'\to\cG\to\cG''\to 0$
of flabby abelian sheaves on $Y$, the resulting sequence
$0\to f_*\cG'\to f_*\cG\to f_*\cG''\to 0$ is still short
exact; the latter follows from \ref{lem_qc-flabby}(i).

(iv): The assertion is local on $X$, so we may assume
that $X$ is quasi-compact and quasi-separated, and then
$Y$ is also quasi-compact. Remark
\ref{rem_acyclic-crit}(iv,v,vi) and lemma
\ref{lem_qc-flabby}(iii.b) reduce to showing that, for
every short exact sequence
$0\to\cG'\to\cG\to\cG''\to 0$ of qc-flabby abelian
sheaves on $Y$, the resulting sequence
$0\to(f_*\cG')_x\to(f_*\cG)_x\to(f_*\cG'')_x\to 0$
is short exact, for every $x\in X$. However, every
such $x$ has a fundamental system of open neighborhoods
consisting of coherent open subsets of $X$, so we come
down to checking that the induced sequence
$0\to\cG'(f^{-1}U)\to\cG(f^{-1}U)\to\cG''(f^{-1}U)\to 0$
is short exact, for every coherent open neighborhood
$U$ of $x$. But under the current assumptions, $f^{-1}U$
is quasi-compact and quasi-separated, so we conclude by
lemma \ref{lem_qc-flabby}(iii.a).

(v): Let $j:U\to X$ be an open immersion, and $\cI$ any
injective $\cA$-module; the counit of adjunction
$j_!j^*\cA\to\cA$ is a monomorphism of
$\cA$-modules, so it induces a surjection
$$
\phi:\Hom_\cA(\cA,\cI)\to\Hom_\cA(j_!j^*\cA,\cI).
$$
We have a natural identification : $\Hom_\cA(\cA,\cI)\isom\cI(X)$;
by adjunction, we have also a natural identification :
$\Hom_\cA(j_!j^*\cA,\cI)\isom\Hom_{j^*\cA}(j^*\cA,j^*\cI)\isom\cI(U)$.
Under these identifications, $\phi$ corresponds to the
restriction map $\cI(X)\to\cI(U)$, whence the contention.

(vi): Let $U\subset X$ be any open subset; we have to
check that the restriction map $\cF(X)\to\cF(U)$ is
surjective. To this aim, fix $s\in\cF(U)$, and denote by
$S$ the set of all pairs $(V,s_V)$ where $V\subset X$
is an open subset containing $U$, and $s_V\in\cF(V)$
is a section whose restriction to $U$ agrees with $s$.
We endow $S$ with a partial ordering, by declaring that
$(V,s_V)\leq(V',s_{V'})$ if $V\subset V'$ and the
restriction of $s_{V'}$ to $V$ agrees with $s_V$.
A standard application of Zorn's lemma shows that
$S$ admits maximal elements, hence let $(V,s_V)$ be
such a maximal element; it suffices to show that
$V=X$. However, if the latter fails, there exists
a noetherian open subset $U'\subset X$ which is not
contained in $V$, and $U'\cap V$ is quasi-compact
(remark \ref{rem_sorite-spectral}(iv)); since $\cF$
is qc-flabby, we may then find $s_{U'}\in\cF(U')$
whose restriction to $U'\cap V$ agrees with the
restriction of $s_V$. Then there exists a unique
$s_{U'\cup V}\in\cF(U'\cap V)$ whose restriction to
$U'$ and $V$ agrees with $s_{U'}$ and respectively
$s_V$; the pair $(U'\cup V,s_{U'\cup V})$ lies in
$S$, contradicting the maximality of $(V,s_V)$.
\end{proof}

\sset\subsubsection{}\label{subsec_lambdasmus}
Let $\Lambda$ be a small cofiltered category; we consider
the datum consisting of:
\begin{itemize}
\item
Two systems
$X_\bullet:=(X_\lambda~|~\lambda\in\Ob(\Lambda))$ and
$Y_\bullet:=(Y_\lambda~|~\lambda\in\Ob(\Lambda))$ of
locally spectral topological spaces, with quasi-compact
and quasi-separated transition maps :
\set\begin{equation}\label{eq_stupid-notate}
\phi_u:X_\mu\to X_\lambda \qquad \psi_u:Y_\mu\to Y_\lambda
\qquad
\text{for every morphism $u:\mu\to\lambda$ in $\Lambda$.}
\end{equation}
\item
A system
$\cF_\bullet:=(\cF_\lambda~|~\lambda\in\Ob(\Lambda))$
where $\cF_\lambda$ is a sheaf on $Y_\lambda$ for every
$\lambda\in\Ob(\Lambda)$, with transition maps :
$$
\psi^{-1}_u\cF_\lambda\to\cF_\mu
\qquad
\text{for every morphism $u:\mu\to\lambda$ in $\Lambda$}.
$$
\item
A system $g_\bullet:=(g_\lambda:Y_\lambda\to
X_\lambda~|~\lambda\in\Ob(\Lambda))$ of quasi-compact and
quasi-separated continuous maps, such that
$g_\lambda\circ\psi_u=\phi_u\circ g_\mu$ for every morphism
$u:\mu\to\lambda$.
\end{itemize}
The (inverse) limits of the systems $X_\bullet$,
$Y_\bullet$ and $g_\bullet$ are representable
by topological spaces and continuous maps
$$
X:=\lim_\Lambda X_\bullet
\qquad
Y:=\lim_\Lambda Y_\bullet
\qquad
g:=\lim_\Lambda g_\bullet:Y\to X
$$
and the induced continuous maps
$$
\phi_\lambda:X\to X_\lambda
\qquad
\psi_\lambda:Y\to Y_\lambda.
$$
are spectral for every $\lambda\in\Lambda$. Moreover,
$\cF_\bullet$ induces a sheaf on $Y$ :
$$
\cF:=
\colim_{\lambda\in\Ob(\Lambda^o)}\psi^*_\lambda\cF_\lambda.
$$

\begin{proposition}\label{prop_dir-im-and-colim}
In the situation of \eqref{subsec_lambdasmus}, the
following holds :
\begin{enumerate}
\item
Suppose that $X_\lambda$ is spectral for every
$\lambda\in\Ob(\Lambda)$. Then the natural map
$$
\colim_{\lambda\in\Ob(\Lambda^o)}\Gamma(Y_\lambda,\cF_\lambda)\to
\Gamma(Y,\cF)
$$
is a bijection.
\item
If $\cF_\bullet$ is a system of abelian sheaves,
then the natural morphisms :
$$
\colim_{\lambda\in\Ob(\Lambda^o)}
\phi^*_\lambda R^ig_{\lambda*}\cF_\lambda\to R^ig_*\cF
$$
are isomorphisms for every $i\in\N$.
\end{enumerate}
\end{proposition}
\begin{proof} By virtue of proposition
\ref{prop_filter-Deligne}, we may assume that $\Lambda$
is a cofiltered partially ordered set, and we let
$(\Lambda^o,\leq)$ be the opposite ordering on $\Lambda$
(see example \ref{ex_universe}(iii)). We first check
the injectivity of the map in (i). Since the system
$Y_\bullet$ is cofiltered, we come down to proving the
following assertion. Let $\lambda\in\Lambda^o$ be any
index, $s,s'\in\Gamma(Y_\lambda,\cF_\lambda)$ two sections
whose images agree in $\Gamma(Y,\cF)$,
and for every $\mu\in\Lambda^o$ such that $\mu\geq\lambda$,
denote by $s_\mu$ and $s'_\mu$ the images of $s$ and
$s'$ in $\Gamma(Y_\mu,\cF_\mu)$; then there exists an
index $\mu\geq\lambda$, such that $s_\mu=s'_\mu$.
However, set
$$
Z_\mu:=\{y\in Y_\mu~|~s_{\mu,y}\neq s'_{\mu,y}\}
\qquad
\text{for every $\mu\in\Lambda^o$ such that $\mu\geq\lambda$}
$$
and we endow $Z_\mu$ with the topology induced from
$Y^c_\mu$ (notation of definition \ref{subsec_constr-top}).
Now $Z_\mu$ is a closed subset of $Y_\mu$ for
every $\mu\geq\lambda$, so it is also closed in
$Y^c_\mu$, and it is easily seen that
$\psi_u(Z_\mu)\subset Z_\nu$ for every morphism
$u:\mu\geq\nu$. Thus, $(Z_\mu~|~\mu\geq\lambda)$
is a cofiltered system of compact topological
spaces (theorem \ref{th_main-spectral}(i)),
and the assumption implies that the limit of this
system is empty. We deduce that $Z_\mu=\emptyset$
for some $\mu\geq\lambda$, so $s_\mu=s'_\mu$ for
such index $\mu$.

Next, we verify the surjectivity. In case
$Y=\emptyset$, we know that $Y_\lambda=\emptyset$
for some $\lambda\in\Lambda^o$, so the assertion holds.
Hence, we may assume that $Y$ is not empty. Now,
let $s\in\Gamma(Y,\cF)$ be any section; since
$$
\cF_y=\colim_{\lambda\in\Lambda^o}\cF_{\lambda,\psi_\lambda(y)}
\qquad
\text{for every $y\in Y$}
$$
we may find, for every $y\in Y$, an index
$\lambda(y)\in\Lambda^o$, a quasi-compact open neighborhood
$U_y$ of $\psi_{\lambda(y)}(y)$ in $Y_{\lambda,y}$ and a
section $s(y)\in\Gamma(U_y,\cF_\lambda)$ such that
the image of $s(y)$ in $\cF_y$ equals $s_y$.
Then, for each $y\in Y$, we may find a quasi-compact
open neighborhood $U'_y$ of $y$ in $\psi^{-1}_{\lambda(y)}U_y$
such that the image of $s(y)$ in $U'_y$ agrees with
the restriction of $s$. The family $(U'_y~|~y\in Y)$
is an open covering of $Y$, and since $Y$
is quasi-compact (theorem \ref{th_main-spectral}(i)), we
may find a finite (non-empty) subset $S\subset Y$
such that $(U'_y~|~y\in S)$ already covers $Y$.
Then, for each $y\in S$ we may find an index
$\mu(y)\geq\lambda(y)$ and a quasi-compact open subset
$U''_y$ of $\psi_{\mu(y)\lambda(y)}^{-1}U_{\lambda(y)}$
such that $\psi^{-1}_{\mu(y)}U''_y=U'_y$
(corollary \ref{cor_main-spectral}(ii.a)).
Next, since $\Lambda$ is cofiltered, we may find
$\mu\in\Lambda^o$ such that $\mu\geq\mu(y)$ for every
$y\in S$, and after replacing $U''_y$ by
$\psi^{-1}_{\mu,\mu(y)}U''_y$, and $s(y)$ by its
image in $\Gamma(\psi^{-1}_{\mu,\mu(y)}U''_y,\cF_\mu)$
for every $y\in S$, we may assume that $\mu(y)$ is a
single index $\mu$, independent of $y\in S$. Now,
since $\bigcup_{y\in S}\psi^{-1}_\mu U''_y=Y$,
there exists $\nu\geq\mu$ such that
$\bigcup_{y\in S}\psi^{-1}_{\nu\mu}U''_y=Y_\mu$
(corollary \ref{cor_main-spectral}(ii.b)), so
we may replace $\mu$ by $\nu$, each $U''_y$
by its preimage in $Y_\nu$, and $s(y)$ by its
image in $\Gamma(\psi^{-1}_{\nu\mu}U''_y,\cF_\nu)$,
and further assume that $(U''_y~|~y\in S)$ is a
covering of $Y_\mu$. Lastly, notice that
-- by construction -- for every $y,y'\in S$,
the images of $s(y)$ and $s(y')$ agree in
$\Gamma(\psi^{-1}_\mu(U_y\cap U_{y'}),\cF)$.
Since the injectivity of the map in (i) has already
been established, we deduce that, for every $y,y'\in S$
there exists $\lambda(y,y')\geq\mu$ such that
the images of $s(y)$ and $s(y')$ agree in
$\Gamma(\psi^{-1}_{\lambda(y,y'),\mu}(U''_y\cap U''_{y'}),
\cF_{\lambda(y,y')})$. Since $\Lambda$ is cofiltered,
we may then find $\lambda\in\Lambda^o$ such that
$\lambda\geq\lambda(y,y')$ for every $y,y'\in S$.
Now, let $V_y:=\psi^{-1}_{\lambda\mu}U''_y$ and
denote by $s'(y)$ the image of $s(y)$ in
$\Gamma(V_y,\cF_\lambda)$, for every $y\in S$; by
construction we have
$s'(y)_{|V_y\cap V_{y'}}=s'(y)_{|V_y\cap V_{y'}}$ for
every $y,y'\in S$, so the system $(s'(y)~|~y\in S)$
corresponds to a well defined section
$s'\in\Gamma(Y_\lambda,\cF_\lambda)$, whose image
in $\Gamma(Y,\cF)$ is finally the original section $s$.

(ii): We notice :

\begin{claim}\label{cl_old-flame}
Suppose that $\cF_\lambda$ is qc-flabby for every
$\lambda\in\Lambda^o$. Then the same holds for $\cF$.
\end{claim}
\begin{pfclaim} Indeed, let $U\subset V$ be any two
quasi-compact open subsets of $Y$; we need to show that the
restriction map
$$
\Gamma(V,\cF)\to\Gamma(U,\cF)
$$
is onto. By corollary \ref{cor_main-spectral}(ii) we can
find $\lambda\in\Lambda^o$ and quasi-compact open subsets
$U_\lambda\subset V_\lambda\subset Y_\lambda$ such that
$U=\psi^{-1}_\lambda U_\lambda$, and likewise for $V$.
Let us set $U_\mu:=\psi^{-1}_{\mu\lambda}U_\lambda$
for every $\mu\in\Lambda^o$ such that $\mu\geq\lambda$,
and likewise define $V_\mu$. Then, up to replacing
$\Lambda^o$ by the cofinal subset
$\{\mu\in\Lambda^o~|~\mu\geq\lambda\}$, we may assume that
$U_\mu\subset V_\mu$ are defined for every $\mu\in\Lambda^o$.
By (i) the natural map
$$
\colim_{\lambda\in\Lambda^o}\Gamma(U_\lambda,\cF_\lambda)\to
\Gamma(U,\cF)
$$
is bijective for every $n\in\N$, and likewise for $V$,
whence the claim.
\end{pfclaim}

Now, we pick for each $\lambda\in\Lambda^o$ an
injective resolution $\cF_\lambda\to\cI_\lambda^\bullet$,
having care to construct $\cI_\lambda^\bullet$ functorially,
so that $\cF_\bullet$ extends to a compatible system
of complexes
$(\cF_\lambda\to\cI_\lambda^\bullet~|~\lambda\in\Lambda^o)$.
From claim \ref{cl_old-flame}, we know that
$\cI^\bullet:=
\colim_{\lambda\in\Lambda^o}\psi^*_\lambda\cI_\lambda^\bullet$
is a complex of qc-flabby sheaves; then lemma
\ref{lem_franziska}(iv) and remark \ref{rem_acyclic-crit}(iv)
yield a natural isomorphism :
$$
g_*\cI^\bullet\isom Rg_*\cF
\qquad
\text{in $\sD(\Z_X\Mod)$}.
$$
To conclude, let $V\to X$ be any quasi-compact
open immersion, which as usual we see as the limit of a system
$(V_\lambda\to X_\lambda~|~\lambda\in\Lambda)$ of quasi-compact
open immersions; we come down to checking that the natural map
$$
\Gamma(V,\colim_{\lambda\in\Lambda^o}
\phi^*_\lambda g_{\lambda*}\cI^\bullet_\lambda)\to
\Gamma(V,g_{\infty*}\cI^\bullet)=\Gamma(g^{-1}V,\cI^\bullet)
$$
is an isomorphism of complexes. However, $g^{-1}V$ is
the limit of the system
$(g^{-1}_\lambda V_\lambda~|~\lambda\in\Lambda)$
of spectral topological spaces and quasi-compact
quasi-separated maps, so (i) naturally identifies
both complexes with
$$
\colim_{\lambda\in\Lambda^o}
\Gamma(g_\lambda^{-1}V_\lambda,\cI^\bullet_\lambda)
$$
whence the contention.
\end{proof}

\begin{remark} In the situation of \eqref{subsec_lambdasmus},
suppose that $X_\bullet$ and $Y_\bullet$ are two systems of
schemes, that the systems $\phi_\bullet$ and $\psi_\bullet$
of \eqref{eq_stupid-notate} consist of affine morphisms of
schemes, and $g_\bullet$ is a system of quasi-compact and
quasi-separated morphisms of schemes.
By \cite[Ch.IV, Prop.8.2.3]{EGAIV-3} the limits (in the
category of schemes) of the systems $X_\bullet$, $Y_\bullet$
(resp. $g_\bullet$) are representable by schemes (resp.
by a morphism of schemes) whose underlying topological
spaces (resp. continuous map) agrees with $X$ and $Y$
(resp. $g$), and the resulting $\phi_\lambda$,
$\psi_\lambda$ are affine morphisms of schemes, for every
$\lambda\in\Ob(\Lambda)$. Hence,  proposition
\ref{prop_dir-im-and-colim} applies especially to such
inverse systems.
\end{remark}

\sset\subsubsection{}\label{subsec_def-Hom-cplx}
Let $(X,\cA)$ be a ringed space, and denote by $\Z_X$
the sheaf of rings associated with the constant presheaf
with values $\Z$; then $\Z_X\Mod$ is the category of
abelian sheaves on $X$, and the natural morphism
$\Z_X\to\cA$ induces a forgetful functor
$$
\phi_\cA:\sD(\cA\Mod)\to\sD(\Z_X\Mod).
$$
Furthermore, let $f:(Y,\cB)\to(X,\cA)$ be a morphism of
ringed spaces. For every $\cB$-module $\cF$, the direct
image $f_*\cF$ is naturally an $\cA$-module, and we notice
that the derived functor
$Rf_*:\sD^+(\cB\Mod)\to\sD^+(\cA\Mod)$ commutes with the
forgetful functors $\phi_\cA,\phi_\cB$; indeed, every
injective $\cB$-module is a flabby $\Z_Y$-module (lemma
\ref{lem_franziska}(v)), and the latter are acyclic for the
functor $f_*:\Z_Y\Mod\to\Z_X\Mod$ (lemma \ref{lem_franziska}(ii)),
so the assertion follows easily.

One defines as in example \ref{ex_monoidal}(v) a
{\em total $\Hom$ cochain complex}, which is a functor :
$$
\cHom_\cA^\bullet:
\sC(\cA\Mod)\times\sC(\cA\Mod)^o\to\sC(\cA\Mod)
$$
on the category of complexes of $\cA$-modules. Recall that,
for any two complexes $M_\bullet$, $N^\bullet$, and any object
$U$ of $X$, the group of $n$-cocycles in
$\cHom^\bullet_{\cA}(M_\bullet,N^\bullet)(U)$ is naturally
isomorphic to $\Hom_{\sC(\cA_{|U}\Mod)}(M_{\bullet|U},N^\bullet_{|U}[n])$,
and $H^n\cHom^\bullet_{\cA}(M_\bullet,N^\bullet)(U)$ is naturally
isomorphic to the group of homotopy classes of maps
$M_{\bullet|U}\to N^\bullet_{|U}[n]$ (see example
\ref{ex_Hom-complex-hots}(i)). We also set :
\set\begin{equation}\label{eq_glob-RHom}
\Hom^\bullet_\cA:=\Gamma\circ\cHom^\bullet_\cA:
\sC(\cA\Mod)\times\sC(\cA\Mod)^o\to\sC(\Gamma(\cA)\Mod).
\end{equation}
The bifunctor $\cHom^\bullet_\cA$ admits a right derived functor :
$$
R\cHom^\bullet_\cA:\sD^+(\cA\Mod)\times\sD(\cA\Mod)^o\to
\sD(\cA\Mod)
$$
for whose construction we refer to \cite[\S10.7]{We}. Likewise,
one has a derived functor $R\Hom_\cA^\bullet$ for \eqref{eq_glob-RHom},
and there are natural isomorphisms of $\cA$-modules :
\set\begin{equation}\label{eq_calculate-D.Hom}
H^iR\Hom^\bullet_\cA(M^\bullet,N^\bullet)\isom
\Hom_{\sD^+(\cA\Mod)}(M^\bullet,N^\bullet[i])
\end{equation}
for every $i\in\Z$, and every bounded below complexes $M^\bullet$
and $N^\bullet$ of $\cA$-modules.

\begin{lemma}\label{lem_triv-dual}
Let $(X,\cA)$ be a ringed space, $\cB$ an $\cA$-algebra.
Then :
\begin{enumerate}
\item
If $\cI$ is an injective $\cA$-module, $\cHom_\cA(\cF,\cI)$
is flabby for every $\cA$-module $\cF$, and $\cHom_\cA(\cB,\cI)$
is an injective $\cB$-module.
\item
There is a natural isomorphism of bifunctors :
$$
R\Gamma\circ R\cHom^\bullet_\cA\isom R\Hom^\bullet_\cA.
$$
\item
The forgetful functor $\sD^+(\cB\Mod)\to\sD^+(\cA\Mod)$
admits the right adjoint :
$$
\sD^+(\cA\Mod)\to\sD^+(\cB\Mod)\quad :\quad
K^\bullet\mapsto R\cHom^\bullet_\cA(\cB,K^\bullet).
$$
\end{enumerate}
\end{lemma}
\begin{proof}(i) : Let $j:U\subset X$ be any open immersion;
a section $s\in\cHom_\cA(\cF,\cI)(U)$ is a map of
$\cA_{|U}$-modules $s:\cF_{|U}\to\cI_{|U}$. We deduce a
map of $\cA$-modules
$j_!s:j_!\cF_{|U}\to j_!\cI_{|U}\to\cI$; since $\cI$
is injective, $j_!s$ extends to a map $\cF\to\cI$, as
required. Next, let $\cB$ be an $\cA$-algebra; we recall
the following :
\begin{claim}\label{cl_right-to-forget} The functor
$$
\cA\Mod\to\cB\Mod\quad :\quad \cF\mapsto\cHom_\cA(\cB,\cF)
$$
is right adjoint to the forgetful functor.
\end{claim}
\begin{pfclaim} We have to exhibit a natural bijection
$$
\Hom_\cA(N,M)\isom\Hom_\cB(N,\Hom_\cA(\cB,M))
$$
for every $\cA$-module $M$ and $\cB$-module $N$. This is given by
the following rule. To an $\cA$-linear map $t:N\to M$ one assigns
the $\cB$-linear map $t':N\to\Hom_\cA(\cB,M)$ such that
$t'(x)(b)=t(x\cdot b)$ for every $x\in N(U)$ and $b\in\cB(V)$,
where $V\subset U$ are any two open subsets of $X$. To describe
the inverse of this transformation, it suffices to remark that
$t(x)=t'(x)(1)$ for every local section $x$ of $N$.
\end{pfclaim}

Since the forgetful functor is exact, the second assertion
of (i) follows immediately from claim \ref{cl_right-to-forget}.
Assertion (ii) follows from (i) and \cite[Th.10.8.2]{We}.

(iii): Let $K^\bullet$ (resp. $L^\bullet$) be a bounded below
complex of $\cB$-modules (resp. of injective $\cA$-modules);
then :
$$
\begin{aligned}
\Hom_{\sD^+(\cA\Mod)}(K^\bullet,L^\bullet)
=\: & H^0\Hom^\bullet_{\sC(\cA\Mod)}(K^\bullet,L^\bullet) \\
=\: & H^0\Hom^\bullet_{\sC(\cB\Mod)}(K^\bullet,\cHom_\cA(\cB,L^\bullet)) \\
=\: & \Hom_{\sD^+(\cB\Mod)}(K^\bullet,\cHom_\cA(\cB,L^\bullet))
\end{aligned}
$$
where the last identity follows from (i) and \eqref{eq_calculate-D.Hom}.
\end{proof}

\begin{theorem}[Trivial duality]\label{th_trivial-dual}
Let $f:(Y,\cB)\to(X,\cA)$ be a morphism of ringed spaces.
For every $M_\bullet\in\Ob(\sD^-(\cA\Mod))$ and
$N^\bullet\in\Ob(\sD^+(\cB\Mod))$, there are natural isomorphisms :
\begin{enumerate}
\item
$R\Hom^\bullet_\cA(M_\bullet,Rf_*N^\bullet)\isom
R\Hom^\bullet_\cB(Lf^*M_\bullet,N^\bullet)$\ \  in\ \  $\sD^+(\cA(X)\Mod)$.
\item
$R\cHom^\bullet_\cA(M_\bullet,Rf_*N^\bullet)\isom
Rf_*R\cHom^\bullet_\cB(Lf^*M_\bullet,N^\bullet)$\ \  in\ \  $\sD^+(\cA\Mod)$.
\end{enumerate}
\end{theorem}
\begin{proof} One applies lemma \ref{lem_triv-dual}(ii) to deduce
(i) from (ii).

(ii) : By standard adjunctions, for any
two complexes $M_\bullet$ and $N^\bullet$ as in the proposition,
we have a natural isomorphism of total $\Hom$ cochain complexes :
\set\begin{equation}\label{eq_std-adjs}
\cHom^\bullet_\cA(M_\bullet,f_*N^\bullet)\isom
f_*\cHom^\bullet_\cB(f^*M_\bullet,N^\bullet).
\end{equation}
Now, let us choose a flat resolution $P_\bullet\isom M_\bullet$
of $\cA$-modules, and an injective resolution $N^\bullet\isom I^\bullet$
of $\cB$-modules. In view of \eqref{eq_std-adjs} and lemma
\ref{lem_triv-dual}(i) we have natural isomorphisms :
$$
Rf_*R\cHom^\bullet_\cB(Lf^*M_\bullet,N^\bullet)\isom
f_*\cHom^\bullet_\cB(f^*P_\bullet,I^\bullet)\stackrel{\sim}{\leftarrow}
\cHom^\bullet_\cA(P_\bullet,f_*I^\bullet)
$$
in $\sD^+(\cA\Mod)$. It remains to show that the
natural map
\set\begin{equation}\label{eq_comp-Homs}
\cHom^\bullet_\cA(P_\bullet,f_*I^\bullet)\to
R\cHom^\bullet_\cA(P_\bullet,f_*I^\bullet)
\end{equation}
is an isomorphism. However, we have two spectral sequences :
$$
\begin{aligned}
E^1_{pq}:=H_p\cHom^\bullet_\cA(P_q,f_*I^\bullet) \Rightarrow\: &
H_{p+q}\cHom^\bullet_\cA(P_\bullet,f_*I^\bullet) \\
F^1_{pq}:=H_pR\cHom^\bullet_\cA(P_q,f_*I^\bullet) \Rightarrow\: &
H_{p+q}R\cHom^\bullet_\cA(P_\bullet,f_*I^\bullet)
\end{aligned}
$$
and \eqref{eq_comp-Homs} induces a natural map of
spectral sequences :
\set\begin{equation}\label{eq_from-E-to-F}
E^1_{pq}\to F^1_{pq}
\end{equation}
Consequently, it suffices to show that \eqref{eq_from-E-to-F}
is an isomorphism for every $p,q\in\N$, and therefore we may
assume that $P_\bullet$ consists of a single flat $\cA$-module
placed in degree zero. A similar argument reduces to the case
where $I_\bullet$ consists of a single injective $\cB$-module
sitting in degree zero. To conclude, it suffices to show :

\begin{claim} Let $P$ be a flat $\cA$-module, $I$ an injective
$\cB$-module. Then the natural map :
$$
\cHom_\cA(P,f_*I)[0]\to R\cHom^\bullet_\cA(P,f_*I)
$$
is an isomorphism.
\end{claim}
\begin{pfclaim}[] $R^i\cHom^\bullet_\cA(P,f_*I)$ is the
sheaf associated with the presheaf on $X$ :
$$
U\mapsto R^i\Hom^\bullet_{\cA_{|U}}(P_{|U},f_*I_{|U}).
$$
Since $I_{|U}$ is an injective $\cB_{|U}$-module,
it suffices therefore to show that $R^i\Hom^\bullet_\cA(P,f_*I)=0$
for $i>0$. However, recall that there is a natural isomorphism :
$$
R^i\Hom^\bullet_\cA(P,f_*I)\simeq\Hom_{\sD(\cA\Mod)}(P,f_*I[-i]).
$$
Since the homotopy category $\Hot(\cA\Mod)$ admits a
left calculus of fractions, we deduce a natural isomorphism :
\set\begin{equation}\label{eq_calc-of-fracts}
R^i\Hom^\bullet_\cA(P,f_*I)\simeq
\colim_{E_\bullet\to P}\rHot_{\cA\Mod}(E_\bullet,f_*I[-i])
\end{equation}
where the colimit ranges over the family of quasi-isomorphisms
$E_\bullet\to P$. We may furthermore restrict the colimit
in \eqref{eq_calc-of-fracts} to the subfamily of all such
$E_\bullet\to P$ where $E_\bullet$ is a bounded above complex
of flat $\cA$-modules, since this subfamily is cofinal.
Every such map $E_\bullet\to P$ induces a commutative diagram :
$$
\xymatrix{
\Hom_\cA(P,f_*I[-i]) \ar[r] \ar[d] &
\rHot_{\cA\Mod}(E_\bullet,f_*I[-i]) \ar[d] \\
\Hom_\cB(f^*P,I[-i]) \ar[r] &
\rHot_{\cB\Mod}(f^*E_\bullet,I[-i])
}
$$
whose vertical arrows are isomorphisms. Since $E_\bullet$
and $P$ are complexes of flat $\cA$-modules, the induced
map $f^*E_\bullet\to f^*P$ is again a quasi-isomorphism;
therefore, since $I$ is injective, the bottom horizontal
arrow is an isomorphism, hence the same holds for the top
horizontal one, and the claim follows easily.
\end{pfclaim}
\end{proof}

\sset\subsubsection{}\label{subsec_filtered-obj}
Suppose that $\cF$ is an abelian group object in a given topos
$X$, and let $(\Fil^n\cF~|~n\in\N)$ be a descending filtration
by abelian subobjects of $\cF$, such that $\Fil^0\cF=\cF$. Set
$\cF^n:=\cF/\Fil^{n+1}\cF$ for every $n\in\N$, and denote
by $\gr^\bullet\cF$ the graded abelian object associated
with the filtered object $\Fil^\bullet\cF$; then we have :

\begin{lemma}\label{lem_filtered-coh}
In the situation of \eqref{subsec_filtered-obj}, fix $q\in\N$
and suppose that the natural map :
$$
\cF\to R\lim_{k\in\N}\cF^k
$$
is an isomorphism in $\sD(\Z_X\Mod)$. Then, the following
conditions are equivalent :
\begin{enumerate}
\alphaenu
\item
$H^q(X,\Fil^n\cF)=0$ for every $n\in\N$.
\item
The natural map $H^q(X,\Fil^n\cF)\to H^q(X,\gr^n\cF)$ vanishes
for every $n\in\N$.
\end{enumerate}
\end{lemma}
\begin{proof} Obviously we have only to show that (b) implies
(a). Moreover, let us set :
$$
(\Fil^n\cF)^k:=\Fil^n\cF/\Fil^{k+1}\cF
\qquad\text{for every $n,k\in\N$ with $k\geq n$}.
$$
There follows, for every $n\in\N$, a short exact sequence of
inverse systems of sheaves :
$$
0\to((\Fil^n\cF)^k~|~k\geq n)\to(\cF^k~|~k\geq n)\to\cF^{n-1}\to 0
$$
(where the right-most term is the constant inverse system
which equals $\cF^{n-1}$ in all degrees, with transition
maps given by the identity morphisms); whence a distinguished
triangle :
$$
R\lim_{k\in\N}\,(\Fil^n\cF)^k\to R\lim_{k\in\N}\cF^k\to\cF^{n-1}[0]\to
R\lim_{k\in\N}\,(\Fil^n\cF)^k[1]
$$
which, together with our assumption on $\cF$, easily implies
that the natural map $\Fil^n\cF\to R\lim_{k\in\N}(\Fil^n\cF)^k$
is an isomorphism in $\sD(\Z_X\Mod)$. Summing up, we may replace
the datum $(\cF,\Fil^\bullet\cF)$ by $(\Fil^n\cF,\Fil^{\bullet+n}\cF)$,
and reduce to the case where $n=0$.

The inverse system $(\cF^n~|~n\in\N)$ defines an abelian
group object of the topos  $X^\N$ (see \cite[\S7.3.4]{Ga-Ra});
whence a spectral sequence:
\set\begin{equation}\label{eq_spec-s-lim}
E^{pq}_2:=\lim_{n\in\N}{}^{\!p}\, H^q(X,\cF^n)\Rightarrow
H^{p+q}(X^\N,\cF^\bullet)\simeq H^{p+q}(X,R\lim_{n\in\N}\cF^n).
\end{equation}
By \cite[Cor.3.5.4]{We} we have $E^{pq}_2=0$ whenever $p>1$,
and, in view of our assumptions, \eqref{eq_spec-s-lim}
decomposes into short exact sequences :
$$
0\to\lim_{n\in\N}{}^{\!1}\,H^{q-1}(X,\cF^n)
\xrightarrow{\alpha_q} H^q(X,\cF)\xrightarrow{\beta_q}
\lim_{n\in\N}\,H^q(X,\cF^n)\to 0
\qquad\text{for every $q\in\N$}
$$
where $\beta_q$ is induced by the natural system of
maps $(\beta_{q,n}:H^q(X,\cF)\to H^q(X,\cF^n)~|~n\in\N)$.

\begin{claim}\label{cl_about-beta}
If (b) holds, then $\beta_q=0$.
\end{claim}
\begin{pfclaim} Suppose that (b) holds; we consider the long
exact cohomology sequence associated with the short exact
sequence :
$$
0\to\Fil^{n+1}\cF\to\Fil^n\cF\to\gr^n\cF\to 0
$$
to deduce that the natural map
$H^q(X,\Fil^{n+1}\cF)\to H^q(X,\Fil^n\cF)$ is onto for every
$n\in\N$; hence the same holds for the natural map
$H^q(X,\Fil^n\cF)\to H^q(X,\cF)$.
By considering the long exact cohomology sequence attached to
the short exact sequence :
\set\begin{equation}\label{eq_comparing}
0\to\Fil^n\cF\to\cF\to\cF^n\to 0
\end{equation}
we then deduce that $\beta_{q,n}$ vanishes for every $n\in\N$,
which implies the claim.
\end{pfclaim}

\begin{claim}\label{cl_about-alpha} If (b) holds, then the inverse
system $(H^{q-1}(X,\cF^n)~|~n\in\N)$ has surjective transition maps.
\end{claim}
\begin{pfclaim} By considering the long exact cohomology
sequence associated with the short exact sequence :
$0\to\gr^n\cF\to\cF^{n+1}\to\cF^n\to 0$,
we are reduced to showing that the boundary map
$\partial_{q-1}:H^{q-1}(X,\cF^n)\to H^q(X,\gr^n\cF)$ vanishes for
every $n\in\N$. However, comparing with \eqref{eq_comparing},
we find that $\partial_{q-1}$ factors through the natural map
$H^q(X,\Fil^n\cF)\to H^q(X,\gr^n\cF)$, which vanishes if (b) holds.
\end{pfclaim}

The lemma follows from claims \ref{cl_about-beta} and
\ref{cl_about-alpha}, and \cite[Lemma 3.5.3]{We}.
\end{proof}

Our next result generalizes to arbitrary coherent spaces
a classical theorem from \cite{Toh}, where it was stated
only for noetherian spaces. We need first to introduce a
general truncation construction for abelian sheaves, that
shall be useful also in section \ref{sec_Cousin}.

\sset\subsubsection{}\label{subsec_loc-closed-supp}
Let $(X,\cA)$ be a ringed topological space, $\phi:V\to X$ a
locally closed immersion, denote by $\bar V$ the topological
closure of $V$ in $X$, and set $\partial V:=\bar V\!\setminus\!V$,
which is also a closed subset of $X$. Let also $i:\bar V\to X$
and $j:V\to\bar V$ be the induced closed and respectively
open immersions. We set
$$
\phi_!:=i_*\circ j_!:\phi^*\cA\Mod\to\cA\Mod
\qquad
t_V:=\phi_!\circ\phi^*:\cA\Mod\to\cA\Mod
$$
and we call $\phi_!$ (resp. $t_V$) the functor of
{\em extension by zero from $V$ to $X$} (resp. of
{\em truncation outside $V$}). Notice that we have
natural identifications
\set\begin{equation}\label{eq_natural-stalks}
t_V\cF_x\isom\left\{\begin{array}{lll}
                   \cF_x & \qquad & \text{if $x\in V$} \\
                   0     & \qquad & \text{otherwise}.
             \end{array}\right.
\end{equation}

\begin{lemma}\label{lem_upside-down}
In the situation of \eqref{subsec_loc-closed-supp}, set
$W:=X\!\setminus\!V$, denote by $\psi:W\to X$ the
resulting locally closed immersion, and let $\cF$ be an
$\cA$-module such that $\psi^*\cF=0$. Then there exists
a unique natural isomorphism
$$
f_\cF:t_V\cF\isom\cF
$$
such that $f_{\cF,x}$ is the natural identification
\eqref{eq_natural-stalks}, for every $x\in V$.
\end{lemma}
\begin{proof} The uniqueness of $f_\cF$ is clear. For the existence,
set $\cG:=i^*\cF$; it suffices to exhibit natural isomorphisms
$$
\cF\isom i_*\cG
\qquad\text{and}\qquad
j_!\phi^*\cF\isom\cG
$$
which induce the corresponding natural identifications on
the stalks.
So we may assume that $\phi$ is either an open or closed
immersion. In case $\phi$ is open, we obtain a natural map
$f:\phi_!\cF\to\cF$ as follows. Recall that $\phi_!\cF$ is
the sheaf associated with the presheaf $\cH$ such that
$\cH(U)=\cF(U)$ for $U\subset V$, and $\cH(U)=0$ otherwise;
then clearly there is a natural monomorphism $\cH\to\cF$
in the category of presheaves on $X$, and the latter
factors uniquely through a map as sought. By construction,
$f_x$ is the identity endomorphism of $\cF_x$, for every
$x\in V$, and $f_x$ is also trivially an isomorphism for
$x\notin V$, since in this case the stalks of both source
and target of $f_x$ vanish. So $f$ is indeed an isomorphism
with the required properties.

Lastly, suppose $\phi$ is a closed immersion; for any open
subset $T$ of $V$, let $C_T$ be the cofiltered set of all
open subsets $U$ of $X$ such that $U\cap V=T$, and set
$$
\phi^{-1}\cH(T):=\colim_{U\in C^o_T}\cH(U)
\qquad
\text{for every presheaf of abelian groups $\cH$ on $X$}.
$$
The rule $T\mapsto\phi^{-1}\cH(T)$ for every open subset
$T$ of $V$ yields a well defined abelian presheaf on $V$,
functorial in $\cH$ (the resulting functor $\phi^{-1}$ from
abelian presheaves on $X$ to abelian presheaves on $V$ is
left adjoint to the direct image functor $\phi_*$ : cp. remark
\ref{rem_was-cofinal}(i)). Now, suppose that $\cH$ is a
$\cA$-module such that $\cH_{|X\setminus V}=0$; in this
case, it is easily seen that $\phi^{-1}\cH$ is a
$\phi^*\cA$-module, and moreover the transition map
$\cH(U)\to\cH(U')$ is an isomorphism, for every
$U,U'\in C_T$ such that $U'\subset U$. Especially, this
applies with $\cH:=\cF$, so we get a natural isomorphism
$$
f_U:\cF(U)\isom i^{-1}\cF(U\cap V)=i_*i^{-1}\cF(U)
\qquad
\text{for every open subset $U\subset X$}
$$
and the system of such maps $f_U$ is the sought isomorphism
$\cF\isom i_*i^*\cF$.
\end{proof}

\begin{proposition}\label{prop_double-truncate}
Let $V$ and $V'$ be two locally closed subsets of $X$.
The following holds : 
\begin{enumerate}
\item
We have a unique isomorphism of functors :
$$
\tau^{V,V'}:t_V\circ t_{V'}\isom t_{V\cap V'}
$$
which induces, for every $\cA$-module $\cF$ and every
$x\in V\cap V'$ a commutative diagram
\set\begin{equation}\label{eq_focus}
{\diagram
t_V\circ t_{V'}\cF_x \ar[rr]^-{\tau^{V,V'}_x} \ar[d]_\alpha 
& & t_{V\cap V''}\cF_x \ar[d]^\beta \\
t_{V'}\cF_x \ar[rr]^\gamma & & \cF_x
\enddiagram}
\end{equation}
where $\alpha$, $\beta$ and $\gamma$ are the natural
identifications of \eqref{eq_natural-stalks}.
\item
Let $V''\subset X$ be a third locally closed subset.
Then the diagram
$$
\xymatrix{ t_V\circ t_{V'}\circ t_{V''}
\ar[rrr]^-{t_V*\tau^{V',V''}} \ar[d]_{\tau^{V,V'}*t_{V''}} & & &
t_V\circ t_{V'\cap V''} \ar[d]^{\tau^{V,V'\cap V''}} \\
t_{V\cap V'}\circ t_{V''} \ar[rrr]^-{\tau^{V\cap V',V''}} & & &
t_{V\cap V'\cap V''}
}$$
commutes.
\item
Suppose that :
\begin{enumerate}
\item
either $V\subset V'$ and the resulting
immersion $V\to V'$ is closed
\item
or else $V'\subset V$ and the resulting immersion
$V'\to V$ is open.
\end{enumerate}
Then there exists a unique morphism of functors
$$
c^{V,V'}:t_{V'}\to t_V
$$
which induces, for every $\cA$-module $\cF$ and every
$x\in V\cap V'$, a commutative diagram
$$
\xymatrix{ t_{V'}\cF_x \ar[rr]^-{c^{V,V'}_x} \ar[rd] & & t_V\cF_x \ar[ld] \\
& \cF_x}$$
whose downward arrows are the natural identifications
\eqref{eq_natural-stalks}.
\item
In the situation of\/ {\em(iii)}, let $V''$ be a third
locally closed subset of $X$. Then the resulting
diagram
$$
\xymatrix{ t_{V''}\circ t_{V'} \ar[rr]^-{t_{V''}*c^{V,V'}}
\ar[d]_{\tau^{V'',V'}} & & t_{V''}\circ t_V \ar[d]^{\tau^{V'',V}} \\
t_{V''\cap V'} \ar[rr]^-{c^{V''\cap V,V''\cap V'}} & & t_{V''\cap V}
}$$
commutes.
\item
Let $V\subset V'\subset V''$ (resp. $V''\subset V'\subset V$)
be a chain of closed (resp. open) immersions of locally
closed subsets of $X$. Then
$$
c^{V,V'}\circ c^{V',V''}=c^{V,V''}.
$$
\end{enumerate}
\end{proposition}
\begin{proof}(i): The uniqueness is clear. For the existence,
set $W:=V\cap V''$, and let $\phi:W\to X$ and
$\psi:X\!\setminus\!W\to X$ be the locally closed immersions;
by inspecting the definition, we find natural isomorphisms
$$
\phi^*\circ t_V\circ t_{V'}\cF\isom\phi^*\cF
\xleftarrow{\sim}\phi^*t_W\cF
\qquad
\text{for every $\cA$-module $\cF$}
$$
which induce on stalks the isomorphisms
$\phi^*(\gamma\circ\alpha)$ and respectively $\phi^*\beta$.
There follows an isomorphism
$\phi_!f:t_W\circ t_V\circ t_{V'}\cF\isom
t_{V''}\circ t_{V''}\cF$. On the other hand, it is easily
seen that $\psi^*\circ t_V\circ t_{V'}\cF$ and
$\psi^*\circ t_W\cF$ both vanish, so lemma
\ref{lem_upside-down} naturally identifies $\phi_!f$
with an isomorphism $t_V\circ t_{V'}\cF\isom t_W\cF$
with the sought property.

(ii): The commutativity of the diagram can be checked
on stalks, and then it follows easily from the commutativity
of \eqref{eq_focus} : details left to the reader.

(iii): Again, the uniqueness of $c^{V,V'}$ is clear.
For the existence, consider first case (a), and write
$V=V'\cap Z$ for some closed immersion $\phi:Z\to X$.
The unit of the adjunction $(\phi^*,\phi_*)$ yields
a natural morphism $\eta:\cF\to t_Z\cF$, whence a
morphism $t_{V'}\eta:t_{V'}\cF\to t_{V'}\circ t_Z\cF$,
which, by (i) is naturally identified with a morphism
$t_{V'}\cF\to t_V\cF$ fulfilling the stated condition.
A similar argument works in case (b) : we write $V'=V\cap U$
for some open immersion $\psi:U\subset X$ and we consider
the counit $\eps:t_U\cF\to\cF$ of the adjunction
$(\psi_!,\psi^*)$; then (i) naturally identifies
$t_V\eps$ with a morphism with the sought properties.

(iv) and (v) can be checked on the stalks : the details
shall be left to the reader.
\end{proof}

\begin{theorem}\label{th_tohoku-vanish}
Let $X$ be any coherent space, $\cF$ any abelian sheaf
on $X$. Then :
\begin{enumerate}
\item
$H^i(X,\cF)=0$ for every $i>\dim X$.
\item
If $X$ is spectral, we have more precisely
$$
H^i(X,\cF)=0
\qquad
\text{for every $i>d_\cF:=\sup(\dim\bar{\{x\}}~|~x\in\Supp\,\cF)$}.
$$
\end{enumerate}
\end{theorem}
\begin{proof} (Here, for any subset $T\subset X$, we
denote by $\bar T$ the topological closure of $T$ in $X$).

\begin{claim}\label{cl_reorganize}
Let $d\in\N$ be any integer, and suppose that (i) holds
for every spectral space $X$ of dimension $\leq d$. Then
(ii) holds for every spectral space $X$ and every
$\Z_X$-module $\cF$ such that $d_\cF\leq d$.
\end{claim}
\begin{pfclaim}
Indeed, let $\cS$ be the set consisting of all pairs
$(U,s)$, where $U\subset X$ is a quasi-compact open
subset, and $s\in\cF(U)$ is any element. For any
$(U,s)\in\cS$, let $\Supp(s):=\{x\in U~|~s_x\neq 0\}$,
notice that $\Supp(s)$ is a closed subset of $U$, and
denote by $i_s:\Supp(s)\to U$ and $j_U:U\to X$ the
resulting closed and respectively open immersions;
the section $s$ is the same as a morphism of
$\Z_U$-modules $\Z_U\to j_U^*\cF$ which factors
through the epimorphism $\Z_U\to i_{s*}\Z_{\Supp(s)}$.
By adjunction, there follows a unique morphism of
$\Z_X$-modules $\phi_{U,s}:j_{U!}i_{s*}\Z_{\Supp(s)}\to\cF$.
For any finite subset $S\subset\cS$, we let
$$
\phi_S:\bigoplus_{(U,s)\in S}j_{U!}i_{s*}\Z_{\Supp(s)}\to\cF
$$
be the unique map of $\Z_X$-modules whose restriction to
$j_{U!}i_{s*}\Z_{\Supp(s)}$ agrees with $\phi_{U,s}$,
for every $(U,s)\in S$. It is easily seen that $\cF$ is
the colimit of the filtered system of its subsheaves
$\cF_S:=\Img\,\phi_S$, for $S$ ranging over all finite
subsets of $\cS$. By proposition \ref{prop_dir-im-and-colim}(i)
we deduce that $H^i(X,\cF)$ is the filtered colimit of
the system of abelian groups $(H^i(X,\cF_S)~|~S\subset\cS)$;
it then suffices to show the assertion for each subsheaf
$\cF_S$. Thus, we may assume from start that there exists
a finite subset $S\subset\cS$ with $\cF=\cF_S$. In this case,
we set
$$
Z:=\bar{\Supp\,\cF}=\bigcup_{(U,s)\in S}\bar{\Supp(s)}
$$
and we let $i_Z:Z\to X$ be the closed immersion. Clearly
we have a natural isomorphism
$$
H^i(X,\cF)\isom H^i(Z,i^*_Z\cF)
\qquad
\text{for every $i\in\N$}
$$
and $Z$ is still a spectral space, by corollary
\ref{cor_procon-is-spec}. Hence, we may replace $X$ by
$Z$ and $\cF$ by $i^*_Z\cF$, and assume from start that
$Z=X$ (and clearly, we still have $\cF=\cF_S$). In view
of our assumption, it then suffices to show that
$$
\dim Z=
\sup\Bigl(\dim\bar{\{x\}}~|~x\in\bigcup_{(U,s)\in S}\Supp(s)\Bigr).
$$
To this aim, in light of remark \ref{rem_irred-comps}(v),
we are reduced to checking that
$\dim\bar{\Supp(s)}=\sup(\dim\bar{\{x\}}~|~x\in\Supp(s))$
for every $(U,s)\in S$. However, notice that $\Supp(s)$ is
pro-constructible in $X$, hence $\bar{\Supp(s)}$ is the
set of all specializations in $X$ of the points of $\Supp(s)$
(proposition \ref{prop_closed-under-spec}(i)). Now the
contention follows from remark \ref{rem_specialize}(iii).
\end{pfclaim}

(i): Let $S_X$ be the sober space associated with $X$, and
$f_X:X\to S_X$ the unit of adjunction (proposition
\ref{prop_soberify}); by remark \ref{rem_sober-permanence},
the space $S_X$ is spectral; moreover the functor
$$
f_X^*:\Z_{S_X}\Mod\to\Z_X\Mod
$$
is an equivalence of categories, and the induced
map $\Gamma(S_X,\cF)\to\Gamma(X,f_X^*\cF)$ is a
bijection for every $\Z_{S_X}$-module $\cF$. It
follows easily that the induced map
$$
H^i(S_X,\cF)\to H^i(X,f_X^*\cF)
$$
is an isomorphism as well; thus, we may replace $X$
by $S_X$, and assume from start that $X$ is spectral.

We argue by induction on $d:=\dim X$. If $d=0$, $X$
is a boolean space (see example \ref{ex_boolean}), so every
sheaf on $X$ is qc-flabby, and the assertion follows from
lemma \ref{lem_franziska}(iv).
We may then assume that $d>0$ and that (i) is already known
for every spectral space of dimension $<d$. Suppose, by way
of contradiction, that $H^{d+1}(X,\cF)\neq 0$ for some
$\Z_X$-module $\cF$, and pick any element $c\neq 0$ in this
cohomology group. For any closed subset $Z\subset X$, let
$i_Z:Z\to X$ be the corresponding closed immersion, denote
by $\cZ$ the system of all closed subsets $Z\subset X$ such
that the image of $c$ does not vanish in $H^{d+1}(Z,i_Z^*\cF)$,
and endow $\cZ$ with the partial ordering given by inclusion.
Notice that each $Z\in\cZ$ is also a spectral space
(corollary \ref{cor_procon-is-spec}), so we must have
$\dim Z=d$, by inductive assumption.

\begin{claim}\label{cl_minimal-support}
(i)\ \ The partially ordered set $\cZ$ admits minimal elements.

\begin{enumerate}
\addenu
\item
Every minimal element of $\cZ$ is an irreducible closed
subset of $X$.
\end{enumerate}
\end{claim}
\begin{pfclaim}(i): Let $(Z_i~|~i\in I)$ be any totally
ordered non-empty subset of $\cZ$; by Zorn's lemma, it suffices
to check that $\bigcap_{i\in I}Z_i\in\cZ$. However, notice
that the inclusion maps $Z_i\to Z_j$ are spectral for every
$i,j\in I$ with $Z_i\subset Z_j$ (lemma \ref{lem_sorite-construct}(x.e)
and remark \ref{rem_sorite-qcoh-maps}(iv)), so the latter
assertion follows directly from proposition
\ref{prop_dir-im-and-colim}(i).

(ii): Let $Z$ be any minimal element of $\cZ$, and
suppose, by way of contradiction, that $Z=Z_1\cup Z_2$
for two non-empty closed subsets $Z_1,Z_2$ strictly
contained in $Z$; notice that $Z_1$ and $Z_2$ are also
spectral spaces (corollary \ref{cor_procon-is-spec}).
Pick any non-empty quasi-compact open subset
$U\subset Z\setminus Z_1$, and set $Z'_1:=Z\setminus U$,
$Z'_2:=\bar U$. Since $U\subset Z_2$, we have
$Z'_2\subset Z_2$ and $Z_1\subset Z'_1$. Moreover, both
$Z'_1$ (resp. $Z'_2$) is the set of all specializations
of the points of $Z\setminus U$ (resp. of $U$); it
follows that $Z'_3:=Z'_1\cap Z'_2$ has dimension $<d$ :
indeed, any totally ordered sequence $x_1<\cdots<x_k$
of elements of $(Z'_3,\leq)$ (with ordering given by
specialization) can be extended to a totally ordered
sequence of $(Z'_2,\leq)$, by adding a generization
$x_{k+1}$ of $x_k$ with $x_{k+1}\in U$.

On the other hand, consider the sequence
\set\begin{equation}\label{eq_ex-closed-supports}
0\to t_{Z'_3}\cF\to t_{Z'_1}\cF\oplus t_{Z'_2}\cF\to t_Z\cF\to 0.
\end{equation}
whose second and third arrows are given by the
natural transformations provided by proposition
\ref{prop_double-truncate}(iii); by considering
the stalks over the points of $Z$, it is easily
seen that \eqref{eq_ex-closed-supports} is short
exact, whence an exact sequence of abelian groups
\set\begin{equation}\label{eq_consequent}
H^d(Z'_3,i_{Z'_3}^*\cF)\to H^{d+1}(Z,i_Z^*\cF)
\xrightarrow{\ \beta\ }
H^{d+1}(Z'_1,i_{Z'_1}^*\cF)\oplus H^{d+1}(Z'_2,i_{Z'_2}^*\cF).
\end{equation}
However, the first term of \eqref{eq_consequent}
vanishes by inductive assumption, hence $\beta$
is injective; especially the image of $c$ cannot
vanish in both $H^{d+1}(Z'_1,i_{Z'_1}^*\cF)$ and
$H^{d+1}(Z'_2,i_{Z'_2}^*\cF)$, contradicting  the
minimality of $Z$.
\end{pfclaim}

By claim \ref{cl_minimal-support}, we may replace
$X$ by a minimal element of $\cZ$, and assume from
start that $X$ is irreducible. In this case, denote
by $\eta\in X$ the maximal point, and let
$j_\eta:\{\eta\}\to X$ be the inclusion map;
there follows a natural map of $\Z_X$-modules
$f:\cF\to j_{\eta*}j^*_\eta\cF$ such that $f_\eta$
is an isomorphism. Since $j^*_\eta\cF$ is obviously
qc-flabby on $\{\eta\}$, lemma \ref{lem_franziska}(iii,iv)
implies that $H^i(X,j_{\eta*}j^*_\eta\cF)=0$ for every $i>0$.
We deduce an exact sequence of abelian groups
\set\begin{equation}\label{eq_scheiderer}
H^{d+1}(X,\Ker\,f)\to H^{d+1}(X,\cF)\to H^{d+1}(X,\Img\,f)
\end{equation}
as well as an isomorphism of abelian groups
$$
H^d(X,\Coker\,f)\isom H^{d+1}(X,\Img\,f)
$$
and by construction we have $d_{\Ker,f},d_{\Coker\,f}<d$.
From the inductive assumption and claim \ref{cl_reorganize},
it follows that the first and third terms of
\eqref{eq_scheiderer} vanish, so the same holds for
the middle one, and the proof is concluded.
\end{proof}

\begin{remark} A different proof of (a slight improvement
of) theorem \ref{th_tohoku-vanish} is found in \cite{Sch}.
See also proposition \ref{prop_patch-trick}(ii).
\end{remark}

\subsection{\v{C}ech cohomology}
This section is a review of the standard constructions
of resolutions via \v{C}ech complexes, for general topoi.
We also include two classical applications to the cohomology
of quasi-coherent modules on a scheme.

\sset\subsubsection{}\label{subsec_Cech-resolution}
Let $T$ be a topos, and $\fU:=(U_i~|~i\in I)$ a family of
objects of $T$; pick any final object $1_T$ of $T$, and
set as well (notation of example \ref{ex_localization-topos}(iii))
$$
Y:=\coprod_{i\in I}U_i
\qquad
U:=\bigcup_{i\in I}\Img(U_i\to 1_T)
\qquad
Z:=\sC U.
$$
According to example \ref{ex_localization-topos}(i,iii) there
exist natural morphisms of topoi
$$
j_U:T/U\to T
\qquad
\pi:T/Y\to T
\qquad
i:Z\to T.
$$
Moreover, the functor $\pi^*:T\to T/Y$ admits a left adjoint
$\pi_!:T/Y\to T$. Furthermore, let $\Z_T$ (resp. $\Z_X$, for
every object $X$ of $T$) be the constant ring object of $T$
(resp. of $T/X$) associated with the ring $\Z$. The functor
$\pi^*:\Z_T\Mod\to\Z_Y\Mod$ admits a left adjoint
$$
\pi_!:\Z_Y\Mod\to\Z_T\Mod.
$$
Then, according to \eqref{subsec_perp-adj}, the adjoint pair
$(\pi_!,\pi^*)$ determines a cotriple $(\perp^\fU,\eta,\mu)$,
which in turns yields a functor
$$
\Z_T\Mod\to\widehat{s}.\Z_T\Mod
\qquad
\cF\mapsto\perp^\fU_\bullet\cF.
$$
(notation of definition \ref{def_simplicial-cats}(iv)).
This construction can be described explicitly as follows.
We attach to $I$ the simplicial set $I^\bullet$ which, in
every degree $n\in\N$ consists of the set of all mappings
$[n]\to I$, where $[n]:=\{0,\dots,n\}$; {\em i.e.}, this
is the cartesian power $I^{n+1}$ of $I$, whose elements
are all the ordered sequences
$\underline t:=(t_0,\dots,t_n)$ of elements of $I$. Every
morphism $\phi:[n]\to[m]$ in the simplicial category
$\Delta$ induces a map $\phi^*:I^{m+1}\to I^{n+1}$ by the rule :
$$
\underline t\mapsto\underline t\circ\phi:=
(t_{\phi 0},\dots,t_{\phi n})
\qquad
\text{for every $\underline t\in I^{m+1}$}
$$
and it is clear that the rules : $[n]\mapsto I^{n+1}$ and
$\phi\mapsto\phi^*$ yield a well defined functor
$\Delta^o\to\Set$. For every $\underline t\in I^{n+1}$, we
fix as well an object $U_{\underline t}$ of $T$ that represents
the product $U_{t_0}\times\cdots\times U_{t_n}$, and denote
by $j_{\underline t}:T/U_{\underline t}\to T$ the induced morphism.
Notice that every $\phi:[n]\to[m]$ in $\Delta$ and every
$\underline t\in I^{m+1}$ induces a well defined morphism in $T$ :
$$
\nu^{(\underline t)}_\phi:U_{\underline t}\to U_{\phi^*(\underline t)}
$$
namely, the unique morphism whose composition with the
projection on the $i$-th factor of $U_{\phi^*(\underline t)}$
agrees with the projection onto the $\phi(i)$-th factor
of $U_{\underline t}$, for every $i=0,\dots,n$. Then
$\nu^{(\underline t)}_\phi$ in turns induces a morphism
$$
u^{(\underline t)}_\phi:j_{\underline t!}j^*_{\underline t}\cF\to
j_{\phi^*(\underline t)!}j^*_{\phi^*(\underline t)}\cF.
$$
With this notation, we have
$$
\perp_n^\fU\cF:=
\bigoplus_{\underline t\in I^{n+1}}j_{\underline t!}j^*_{\underline t}\cF
\qquad
\text{for every $n\in\N$}
$$
and every morphism $\phi:[n]\to[m]$ in $\Delta$ corresponds to
the unique morphism of $\Z_T$-modules
$u_\phi:\perp^\fU_m\cF\to\perp^\fU_n\cF$ that makes commute the
diagram
$$
{\diagram j_{\underline t!}j^*_{\underline t}\cF
\ar[rr]^-{u^{(\underline t)}_\phi} \ar[d] & &
j_{\phi^*(\underline t)!}j^*_{\phi^*(\underline t)}\cF \ar[d] \\
\perp_m^\fU\cF \ar[rr]^-{u_\phi} & & \perp_n^\fU\cF
\enddiagram}
\qquad
\text{for every $\underline t\in I^{m+1}$}.
$$
where the vertical arrows are the inclusion maps. Moreover,
$\perp^\fU_{-1}\cF=\cF$, and the augmentation
$\eps:\perp_0^\fU\cF\to\perp^\fU_{-1}(\fU)$ is the sum of the
natural morphisms $j_{i!}j^*_i\cF\to\cF$ given by the counit
for the adjunction $(j_{i!},j^*_i)$, for every $i\in I$.

\begin{definition}\label{def_Cech-complex} 
(i)\ \
With the notation of \eqref{subsec_Cech-resolution}, the
{\em augmented \v{C}ech resolution} associated with $\fU$
and the $\Z_T$-module $\cF$ is the unnormalized chain complex
of $\Z_T$-modules
$$
(R_\bullet(\fU,\cF),d_\bullet)
$$
obtained from the augmented simplicial complex
$\perp^\fU_\bullet(\fU,\Z_T)$ (see definition
\ref{def_Dold-Kan}(i)).

(ii)\ \
To ease notation, we shall set also
$$
R_\bullet(\fU):=R_\bullet(\fU,\Z_T)
\qquad\text{and}\qquad
\bar R_\bullet(\fU):=\tau_{\leq 0}R_\bullet(\fU).
$$
Then the
{\em augmented \v{C}ech complex of $\cF$ relative to the
covering $\fU$} of $U$ is the cochain complex
$$
C^\bullet(\fU_\bullet,\cF):=
\Hom^\bullet_{\Z_T}(R_\bullet(\fU),\cF[0]).
$$
Hence, $R_\bullet(\fU)$ is an object of $\sC^{\leq 1}(\Z_T\Mod)$,
and $C^\bullet(\fU_\bullet,\cF)$ is an object of
$\sC^{\geq -1}(\Z\Mod)$.
\end{definition}

\begin{lemma}\label{lem_cotriples-are-back}
With the notation of definition {\em\ref{def_Cech-complex}},
we have :
\begin{enumerate}
\item
For every $\Z_T$-module $\cF$, the complex of\/
$\Z_U$-modules $j_U^*R_\bullet(\fU,\cF)$ is acyclic.
\item
The natural projection $\Z_T\to i_*\Z_Z$ induces a
quasi-isomorphism :
$$
R_\bullet(\fU)\isom i_*\Z_Z[-1]
\qquad
\text{in $\sC^-(\Z_T\Mod)$}.
$$
\item
Suppose that $\fU$ is a {\em covering} of\/ $T$ ({\em i.e.}
that $U=1_T$ and $Z=\emptyset_T$). Then the differential $d_0$
of $R_\bullet(\fU)$ induces a quasi-isomorphism
$$
\bar R_\bullet(\fU)\isom\Z_T[0]
\qquad
\text{in $\sC^-(\Z_T\Mod)$}.
$$
\end{enumerate}
\end{lemma}
\begin{proof}(i): The natural morphism $Y\to U$ is an
epimorphism, hence it is a covering in the canonical
topology of $T$, and therefore it suffices to check that
$\pi^*R_\bullet(\fU,\cF)$ is acyclic. But proposition
\ref{prop_triple-res}(i) says that
$\pi^*(\perp^\fU_\bullet\cF)$ is homotopically trivial. Then
the assertion follows from corollary \ref{cor_Dold-Kan}(ii).

Assertions (ii) and (iii) are immediate consequences of (i).
\end{proof}

\sset\subsubsection{}\label{subsec_choose-total}
Keep the notation of \eqref{subsec_Cech-resolution}, and
choose any total ordering on $I$. We let
$I_\alt^{n+1}\subset I^{n+1}$ be the subset of all injective
order-preserving mappings $[n]\to I$, for every $n\in\N$,
and we set
$$
R^\alt_{-1}(\fU,\cF):=\cF
\qquad
R_n^\alt(\fU,\cF):=\bigoplus_{\underline t\in I_\alt^{n+1}}
j_{\underline t!}j^*_{\underline t}\cF
\qquad
\text{for every $n\in\N$}.
$$
Obviously $\phi^*(I_\alt^{m+1})\subset I_\alt^{n+1}$ for every
{\em injective} morphism $\phi:[n]\to[m]$ of $\Delta$, so
the differential $d_n$ of $R_\bullet(\fU,\cF)$ restricts to
a well defined morphism of $\Z_T$-modules
$R^\alt_n(\fU,\cF)\to R^\alt_{n-1}(\fU,\cF)$ for every $n\in\N$,
and we obtain therefore a subcomplex
$$
(R^\alt_\bullet(\fU,\cF),d_\bullet).
$$
We set as well $R^\alt_\bullet(\fU):=R^\alt_\bullet(\fU,\Z_T)$, and
to every $\Z_T$-module $\cF$ we assign its {\em alternating
augmented \v{C}ech complex relative to the covering $\fU$}
of $U$, which is the cochain complex
$$
C^\bullet_\mathrm{alt}(\fU_\bullet,\cF):=
\Hom^\bullet_\Z(R^\alt_\bullet(\fU),\cF[0]).
$$
We shall also consider the alternating variant of the
{\em \v{C}ech resolution}
$$
\bar R{}^\alt_\bullet(\fU):=\tau_{\leq 0}R^\alt_\bullet(\fU).
$$

\sset\subsubsection{}\label{subsec_explain-alt}
For every $n\in\N$, denote by $S_{n+1}$ the permutation group
of $[n]$; then, for any $\underline t\in I_\alt^{n+1}$, and any
$\sigma\in S_{n+1}$, the sequence $\underline t\circ\sigma$ is
an element of $I^{n+1}$ such that
$U_{\underline t\circ\sigma}=U_{\underline t}$, and there exists a
unique isomorphism
$\mu_{\underline t}:U_{\underline t\circ\sigma}\isom U_{\underline t}$
whose composition with the projection on the $i$-th factor
of $U_{\underline t}$ agrees with the projection on the
$\sigma(i)$-th factor of $U_{\underline t\circ\sigma}$, for every
$i\in[n]$. We set
$$
\mathrm{Alt}^{(\underline t\circ\sigma)}:=
\sign(\sigma)\cdot v_{\underline t}:
j_{(\underline t\circ\sigma)!}\Z_{U_{\underline t\circ\sigma}}\to
j_{\underline t!}\Z_{U_{\underline t}}
$$
where $v_{\underline t}$ is the isomorphism induced by
$\mu_{\underline t}$, and where $\sign(\sigma)$ denotes
the signature of $\sigma$. Then, for every $n\in\N$
there exists a unique morphism of $\Z_T$-modules
$$
\mathrm{Alt}^\fU_n:R_n(\fU)\to R_n^\alt(\fU)
$$
whose kernel contains the direct summands
$j_{\underline t!}\Z_{U_{\underline t}}$ such that the mapping
$\underline t:[n]\to I$ is not injective, and that makes
commute the diagram
$$
{\diagram
j_{(\underline t\circ\sigma)!}\Z_{U_{\underline t\circ\sigma}}
\ar[rr]^-{\mathrm{Alt}^{(\underline t\circ\sigma)}} \ar[d] & &
j_{\underline t!}\Z_{U_{\underline t}} \ar[d] \\
R_n(\fU) \ar[rr]^-{\mathrm{Alt}^\fU_n} & &R_n^\alt(\fU)
\enddiagram}
\qquad
\text{for every $\underline t\in I^\alt_n$}
$$
(where the vertical arrows are the inclusion maps).

\begin{proposition}\label{prop_Cech-resolve}
In the situation of \eqref{subsec_explain-alt}, denote
by $\iota_\bullet^\fU:R^\alt_\bullet(\fU)\to R_\bullet(\fU)$
the inclusion map, and suppose that $U_i$ is a subobject
of $1_T$, for every $i\in I$. Then the following holds :
\begin{enumerate}
\item
The system $(\mathrm{Alt}^\fU_{n-1}~|~n\in\N)$ defines an
epimorphism of chain complexes
$$
\mathrm{Alt}^\fU_\bullet:R_\bullet(\fU)\to R^\alt_\bullet(\fU)
$$
whose kernel is independent of the choice of ordering on $I$.
\item
$\mathrm{Alt}^\fU_\bullet\circ\iota^\fU_\bullet=\one_{R^\alt_\bullet(\fU)}$.
and $\iota^\fU_\bullet\circ\mathrm{Alt}^\fU_\bullet$ is homotopically
equivalent to $\one_{R_\bullet(\fU)}$.
\item
Especially, $\iota^\fU_\bullet$ and $\mathrm{Alt}^\fU_\bullet$
are mutually inverse isomorphisms in $\sD(\Z_X\Mod)$.
\item
The natural projection $\Z_T\to i_*\Z_Z$ induces a
quasi-isomorphism :
$$
R^\alt_\bullet(\fU)\isom i_*\Z_Z[-1]
\qquad
\text{in $\sC^-(\Z_T\Mod)$}.
$$
\item
Suppose that $\fU$ is a {\em covering} of\/ $T$ ({\em i.e.}
that $U=1_T$ and $Z=\emptyset_T$). Then the differential $d_0$
of $R^\alt_\bullet(\fU)$ induces a quasi-isomorphism
$$
\bar R{}^\alt_\bullet(\fU)\isom\Z_T[0]
\qquad
\text{in $\sC^-(\Z_T\Mod)$}.
$$
\end{enumerate}
\end{proposition}
\begin{proof}(i): Fix an integer $n>0$, a sequence
$\underline s\in I^{n+1}$, and let $c\leq n+1$ be the
cardinality of $\{s_0,\dots,s_n\}$. It suffices to check
that the compositions $d_n\circ\mathrm{Alt}^\fU_n$ and
$\mathrm{Alt}^\fU_{n-1}\circ d_n$ agree on
$j_{\underline s!}\Z_{U_{\underline s}}$. In case $c<n$, it
is easily seen that both these compositions vanish
on such direct summand.
In case $c=1$ and $n=0$, the assertion follows by a
direct inspection.
In case $c=n+1>1$, we may find a unique permutation
$\sigma\in S_{n+1}$ and a unique $\underline t\in I^{n+1}_\alt$
such that $\underline s=\underline t\circ\sigma$, and we have
to show that we have a commutative diagram
$$
\xymatrix{ j_{(\underline t\circ\sigma)!}\Z_{U_{\underline t\circ\sigma}}
\ar[rr]^-{\mathrm{Alt}^{(\underline t\circ\sigma)}}
\ar[d]_{d_{\underline t\circ\sigma}} & &
j_{\underline t!}\Z_{U_{\underline t}} \ar[d]^{d_{\underline t}} \\
R_{n-1}(\fU) \ar[rr]^-{\mathrm{Alt}^\fU_{n-1}} & & R_{n-1}(\fU')
}$$
where $d_{\underline t}$ and $d_{\underline t\circ\sigma}$
denote the restrictions of the differentials of
$R^\alt_\bullet(\fU)$ and respectively $R_\bullet(\fU)$.
To this aim, for every $i=0,\dots,n$ denote by $\beta^{(i)}_1$
the unique permutation of $[n]$ such that
$$
\beta^{(i)}_1(\sigma(i))=\sigma(i)
\qquad\text{and}\qquad
\beta^{(i)}_1(\sigma(j))<\beta^{(i)}_1(\sigma(k))
\quad
\text{for every $j,k\in[n]\setminus\{i\}$ with $j<k$}.
$$
Also, we define a third permutation as follows :
\begin{itemize}
\item
If $\sigma(i)<i$, let $\beta^{(i)}_2$ be the cycle
$(i,i-1,\dots,\sigma(i))$ of length $i-\beta(i)+1$
\item
If $\sigma(i)>i$, let $\beta^{(i)}_2$ be the cycle
$(i,i+1,\dots,\sigma(i))$ of length $\sigma(i)-i+1$
\item
If $\sigma(i)=i$, let $\beta^{(i)}_2$ be the identity map
of $[n]$.
\end{itemize}
It is easily seen that in either cases
$$
\sigma^{-1}=\beta^{(i)}_2\circ\beta^{(i)}_1
\qquad\text{and}\qquad
\sign(\beta^{(i)}_2)=(-1)^{i-\sigma(i)}
\qquad
\text{for every $i=0,\dots,n$}.
$$
Now, for every $i=0,\dots,n$ let $\eps_i:[n-1]\to[n]$
be the $i$-th face map (see \eqref{subsec_do-faces});
we compute:
$$
\begin{aligned}
\mathrm{Alt}^\fU_{n-1}\circ d_{\underline t\circ\sigma}=\,&
\mathrm{Alt}^\fU_{n-1}\circ\Bigl(
\sum_{i=0}^n(-1)^i\cdot u^{(\underline t\circ\sigma)}_{\eps_i}\Bigr) \\
=\,&
\sum_{i=0}^n(-1)^i\cdot\sign(\beta^{(i)}_1)\cdot
u^{(\underline t\circ\sigma)}_{\eps_i} \\
=\,&
\sum_{i=0}^n(-1)^{\sigma(i)}\cdot\sign(\sigma)\cdot
u^{(\underline t\circ\sigma)}_{\eps_i} \\
=\,&
\sum_{i=0}^n(-1)^i\cdot\sign(\sigma)\cdot u^{(\underline t)}_{\eps_i} \\
=\,&
\sign(\sigma)\cdot d_{\underline t}
\end{aligned}
$$
as required. Lastly, if $c=n$, the morphism
$d_n\circ\mathrm{Alt}^{(\underline t\circ\sigma)}$ vanishes, and
it remains to check that
$\mathrm{Alt}^\fU_{n-1}\circ d_{\underline s}=0$ as well. However,
in this situation there exist exactly two distinct indices
$i,j\in[n]$ such that $s_i=s_j$, and it follows already that
$\mathrm{Alt}^\fU_{n-1}\circ u^{(\underline s)}_{\eps_k}=0$ for
every $k\neq i,j$. Moreover, clearly there exists a unique
$\underline t\in I^n_\alt$ such that
$\{t_0,\dots,t_{n-1}\}=\{s_0,\dots,s_n\}$, and notice that
the diagram :
$$
\xymatrix{ U_{\underline s} \ar[rr]^-{\nu^{(\underline s)}_{\eps_i}}
\ar[d]_{\nu^{(\underline s)}_{\eps_j}} & &
U_{\eps^*_i(\underline s)} \ar[d]^{\mu_{\eps^*_i(\underline s)}} \\
U_{\eps^*_j(\underline s)} \ar[rr]^{\mu_{\eps^*_j(\underline s)}} & &
U_{\underline t}
}$$
commutes (this is the only place where we use our assumption
on the objects $U_i$). There follows a commutative diagram
$$
\xymatrix{ j_{\underline s!}\Z_{U_{\underline s}}
\ar[rr]^-{u^{(\underline s)}_{\eps_i}}
\ar[d]_{u^{(\underline s)}_{\eps_j}} & &
j_{\eps^*_i(\underline s)!}\Z_{U_{\eps^*_i(\underline s)}}
\ar[d]^{v_{\eps^*_i(\underline s)}} \\
j_{\eps^*_j(\underline s)!}\Z_{U_{\eps^*_j(\underline s)}}
\ar[rr]^{v_{\eps^*_j(\underline s)}} & &
U_{\underline t}
}$$
Now, let $\beta_i$ (resp. $\beta_j$) be the unique permutation
of $[n]$ that fixes $i$ (resp. $j$) and such that
$\underline s\circ\beta_i$ (resp. $\underline s\circ\beta_j$)
restricts to an order-preserving map on $[n]\setminus\{i\}$
(resp. on $[n]\setminus\{j\}$). Arguing as in the foregoing,
we see that $\beta_i=\gamma\circ\beta_j$ for a cyclic
permutation $\gamma$ of length $|j-i|$, so that
$$
\rho:=(-1)^i\cdot\sign(\beta_i)=(-1)^{j+1}\cdot\sign(\beta_j).
$$
Denote by $\one_{\underline t}$ the identity automorphism of
$j_{\underline t!}\Z_{U_{\underline t}}$. Summing up, we get
$$
\begin{aligned}
\mathrm{Alt}^\fU_{n-1}\circ d_{\underline s}=\,&
(-1)^i\cdot\mathrm{Alt}^\fU_{n-1}\circ u_{\eps_i}^{(\underline s)}+
(-1)^j\cdot\mathrm{Alt}^\fU_{n-1}\circ u_{\eps_j}^{(\underline s)} \\
=\,&
(-1)^i\cdot\sign(\beta_i)\cdot
v_{\eps^*_i(\underline s)}\circ u_{\eps_i}^{(\underline s)}+
(-1)^j\cdot\sign(\beta_j)\cdot
v_{\eps^*_j(\underline s)}\circ u_{\eps_j}^{(\underline s)} \\
=\,& \rho\cdot(v_{\eps^*_i(\underline s)}\circ u_{\eps_i}^{(\underline s)}-
v_{\eps^*_j(\underline s)}\circ u_{\eps_j}^{(\underline s)}) \\
=\,& 0
\end{aligned}
$$
and the proof of (i) is concluded.

(ii): Let $\sV$ be a universe such that $T$ is $\sV$-small,
and set $T':=T^\wedge_\sV$; also, let $J$ be the canonical
topology on $T$. Then $T'$ is a $\sV$-topos, and the Yoneda
imbedding $h:T\to T'$ factors through an equivalence
$T\isom(T,J)^\sim_\sU$ and the inclusion functor
$(T,J)^\sim_\sU\to T'$. Set $U'_i:=h_{U_i}\in\Ob(T')$ for every
$i\in I$, $U':=\bigcup_{i\in I}\Img(U'_i\to 1_{T'})$, and
$\fU':=(U'_i~|~i\in I)$. Hence $U'$ is a presheaf on $T$
whose associated sheaf is isomorphic to $h_U$.
Let $j_{U'}:T'/U'\to T'$ be the induced morphism of topoi;
according to lemma \ref{lem_cotriples-are-back}, the chain
complex $j_{U'}^*R_\bullet(\fU')$ is acyclic, hence the same
holds for $j_{U'!}j_{U'}^*R_\bullet(\fU')$, and there follows
a commutative diagram of complexes of $\Z_{T'}$-modules :
$$
\xymatrix{ \bar R_\bullet(\fU') \ar[r]
\ar[d]_{\tau_{\leq 0}(\iota^{\fU'}_\bullet\circ\mathrm{Alt}^{\fU'}_\bullet)} &
j_{U'!}\Z_{U'}[0] \ddouble \\
\bar R_\bullet(\fU') \ar[r] & j_{U'!}\Z_{U'}[0]
}$$
whose horizontal arrows are induced by the differential
$d_0$ of $R_\bullet(\fU')$, and therefore are isomorphisms
in the category $\sD(\Z_{T'}\Mod)$, so that
$$
\tau_{\leq 0}(\iota^{\fU'}_\bullet\circ\mathrm{Alt}^{\fU'}_\bullet)=
\one_{\bar R_\bullet(\fU')}
\qquad
\text{in $\sD(\Z_{T'}\Mod)$}.
$$
However, from the explicit description of
\eqref{subsec_Cech-resolution} it is easily seen that
$\bar R_\bullet(\fU')$ is a complex of projective
$\Z_{T'}$-modules, hence we have as well
$\tau_{\leq 0}(\iota^{\fU'}_\bullet\circ\mathrm{Alt}^{\fU'}_\bullet)=
\one_{\bar R_\bullet(\fU')}$ in $\Hot(\Z_{T'}\Mod)$, by
theorem \ref{th_enough-is-enough}(iii), and then it follows
as well that
$\iota^{\fU'}_\bullet\circ\mathrm{Alt}^{\fU'}_\bullet=\one_{R_\bullet(\fU')}$
in $\Hot(\Z_{T'}\Mod)$.
Next, notice that $R_\bullet(\fU')$ is a complex
of abelian presheaves on $T$, whose associated complex
of sheaves agrees with $R_\bullet(\fU)$, and likewise
the morphism of abelian sheaves associated with
$\iota^{\fU'}_\bullet\circ\mathrm{Alt}^{\fU'}_\bullet$ agrees
with $\iota^\fU_\bullet\circ\mathrm{Alt}^\fU_\bullet$. The
assertion follows.

(iii) follows directly from (ii). Likewise, (iv) and (v)
are immediate from (iii) and lemma
\ref{lem_cotriples-are-back}(ii,iii).
\end{proof}

\begin{remark}\label{rem_justify-name}
(i)\ \
In light of proposition \ref{prop_Cech-resolve}(iv) we shall
call $R^\alt_\bullet(\fU)$ the {\em alternating augmented \v{C}ech
resolution} associated with $\fU$. This name can be justified
as follows. Notice that for every $\Z_T$-module $\cF$ we have
natural identifications
\set\begin{equation}\label{eq_naturally-cech}
C^n(\fU,\cF)\isom
\prod_{\underline t\in I^{n+1}}\cF(U_{\underline t})
\qquad
C^n_\alt(\fU,\cF)\isom
\prod_{\underline t\in I^{n+1}_\alt}\cF(U_{\underline t})
\end{equation}
for every $n\in\N$ (notation of definition
\ref{def_Cech-complex}(ii)). Under these identifications,
$C^\bullet(\fU,\cF)$ becomes the cochain complex
\set\begin{equation}\label{eq_r-column}
0\to\Gamma(\cF)\to\prod_{t\in I}\cF(U_t)\to\cdots\to
\prod_{\underline t\in I^{n+1}}\cF(U_{\underline t})\to\cdots
\end{equation}
whose differential $d^n$ is given by the following rule.
First, $d^{-1}$ assigns to each global section
$s\in\Gamma(\cF)$ the system of its restrictions
$(s_{|U_t}~|~t\in I)$. In case $n\geq 0$, we have
$$
d^n(s_\bullet)_{\underline t}:=\sum^{n+1}_{k=0}
(-1)^{k+n+1}\cdot(s_{\underline t\circ\eps_k})_{|U_{\underline t}}
\qquad
\text{for every $\underline t\in I^{n+2}$ and every
$s_\bullet\in C^n(\fU,\cF)$}.
$$

(ii)\ \
Likewise, under the natural identifications of (i), the
alternating augmented \v{C}ech complex $C^\bullet_\alt(\fU,\cF)$
becomes the cochain complex
$$
0\to\Gamma(\cF)\to\prod_{t\in I}\cF(U_t)\to\cdots\to
\prod_{\underline t\in I^{n+1}_\alt}\cF(U_{\underline t})\to\cdots
$$
whose differentials are given by the same expressions as in (i).

(iii)\ \
Under the identification of (i), the morphisms
$R_n(\fU)\to\cF$ that factor through $\mathrm{Alt}^\fU_n$
correspond to the system of sections
$f_\bullet:=(f_{\underline t}~|~\underline t\in I^{n+1})$ with
the following properties :
\begin{itemize}
\item
$f_{\underline t}=0$ whenever the map $\underline t:[n]\to I$
is not injective.
\item
$f_{\underline t\circ\sigma}=\sign(\sigma)\cdot f_{\underline t}$
for every injective map $\underline t:[n]\to I$ and every
$\sigma\in S_{n+1}$.
\end{itemize}
Thus, under the identification \eqref{eq_naturally-cech},
the submodule $\Hom_{\Z_X}(R^\alt_n(\fU),\cF)$ corresponds to
the subgroup of all {\em alternating systems $f_\bullet$ of
sections of $\cF$}.
\end{remark}

\sset\subsubsection{}\label{subsec_maddening-combinatorics}
In the situation of \eqref{subsec_Cech-resolution}, consider
another family $\fU':=(U'_{i'}~|~i'\in I')$ of objects of
$T$, indexed by a totally ordered set $I'$, and set
$U':=\bigcup_{i'\in I'}\Img(U'_{i'}\to 1_T)$ and $Z':=\sC U'$.
Also, for every $n\in\N$, and every $\underline s\in I'^{n+1}$
set as usual $U'_{\underline s}:=U'_{s_0}\times\cdots\times U'_{s_n}$,
and denote by $j'_{\underline s}:T/U'_{\underline s}\to T$ the
induced functor. Suppose that $\fU$ is a {\em refinement}
of $\fU'$, {\em i.e.} there exists a mapping (that does
not necessarily respect the orderings)
$$
\tau:I\to I'
\qquad\text{such that}\qquad
\Hom_T(U_i,U'_{\tau(i)})\neq\emptyset
\qquad
\text{for every $i\in I$}
$$
and pick a morphism $\theta_i:U_i\to U'_{\tau(i)}$ for
every $i\in I$. Then, for every $n\in\N$ and every 
$\underline t:=(t_0,\dots,t_n)\in I^{n+1}$, we set
$\tau(\underline t):=(\tau(t_0),\dots,\tau(t_n))\in I'^{n+1}$,
and notice that the morphism
$\theta_{\underline t}:=\theta_{t_0}\times\cdots\times\theta_{t_n}:
U_{\underline t}\to U'_{\tau(\underline t)}$ induces a morphism
of $\Z_T$-modules
$$
\phi_{\underline t}:j_{\underline t!}\Z_{U_{\underline t}}\to
j'_{\tau(\underline t)!}\Z_{U'_{\tau(\underline t)}}.
$$
For every $n\in\N$, denote by $\phi_n:R_n(\fU)\to R_n(\fU')$
the unique morphism of $\Z_T$-modules that makes commute the
diagram
$$
{\diagram
j_{\underline t!}\Z_{U_{\underline t}} \ar[r]^-{\phi_{(\underline t)}}
\ar[d] & j'_{\tau(\underline t)!}\Z_{U'_{\tau(\underline t)}} \ar[d] \\
R_n(\fU) \ar[r]^-{\phi_n} & R_n(\fU')
\enddiagram}
\qquad
\text{for every $\underline t\in I^{n+1}$}.
$$
whose vertical arrows are the inclusion maps. Let as
well $\phi_{-1}:=\one_{\Z_T}$. With this notation, it is
clear that the system $(\phi_{n-1}~|~n\in\N)$ is a morphism
of augmented simplicial complexes
$R_\bullet(\fU)\to R_\bullet(\fU')$, and we denote
$$
\phi_\bullet:R_\bullet(\fU)\to R_\bullet(\fU')
$$
the associated morphism of chain complexes.

\sset\subsubsection{}\label{subsec_alt-refine-map}
Suppose next that $U'_{i'}$ is a subobject of $1_T$,
for every $i'\in I'$; then, on account of proposition
\ref{prop_Cech-resolve}(i) we may define
$$
\phi^\alt_\bullet:=
\mathrm{Alt}^{\fU'}_\bullet\circ\phi_\bullet\circ\iota^\fU_\bullet:
R^\alt_\bullet(\fU)\to R^\alt_\bullet(\fU').
$$
The terms $\phi^\alt_n$ can be described explicitly : first, it
is clear that $\phi^\alt_{-1}:=\one_{\Z_T}$. Next, if $n\in\N$
and $\underline t\in I_\alt^{n+1}$ is any sequence, let
$$
\phi^\alt_{\underline t}:j_{\underline t!}\Z_{U_{\underline t}}\to
R^\alt_n(\fU')
$$
be the map given by the following rule. If
$\tau(\underline t)$ is also an injective sequence
$[n]\to I'$, let $\sigma\in S_{n+1}$ be the unique
permutation such that
$\tau(\underline t)\circ\sigma\in I'^{n+1}_\alt$, and set
$\phi^\alt_{\underline t}:=\sign(\sigma)\cdot\phi_{(\underline t)}:
j_{\underline t!}\Z_{U_{\underline t}}\to\Z_{U'_{\tau(\underline t)\circ\sigma}}\subset
R^\alt_n(\fU')$. In case $\tau(\underline t)$ is not injective,
we let $\phi^\alt_{\underline t}$ be the zero morphism. Then
$\phi^\alt_n$ is the sum of the $\phi^\alt_{\underline t}$,
for $\underline t$ ranging over all elements of $I^{n+1}_\alt$.
In case also $U_i$ is a subobject of $1_T$ for every $i\in I$,
it is easily seen that we get a commutative diagram
$$
\xymatrix{
R_\bullet(\fU) \ar[rr]^-{\mathrm{Alt}^\fU_\bullet} \ar[d]_{\phi_\bullet}
& & R^\alt_\bullet(\fU) \ar[d]^{\phi^\alt_\bullet} \\
R_\bullet(\fU') \ar[rr]^-{\mathrm{Alt}^{\fU'}_\bullet} & &
R^\alt_\bullet(\fU').
}$$

\begin{lemma}\label{lem_mad-combinatorics}
In the situation of  \eqref{subsec_maddening-combinatorics},
suppose that $\tau':I\to I'$ is another mapping with
$\Hom_T(U_i,U'_{\tau'(i)})\neq\emptyset$ for every $i\in I$.
Pick a system of morphisms
$(\theta'_i:U_i\to U'_{\tau'(i)}~|~i\in I)$, and denote by
$\phi'_\bullet:R_\bullet(\fU)\to R_\bullet(\fU')$ the associated
morphisms of chain complexes. We have :
\begin{enumerate}
\item
There exists a homotopy from $\phi_\bullet$ to $\phi'_\bullet$.
\item
If $U'_{i'}$ is a subobject of $1_T$ for every $i'\in I'$,
there is also a homotopy from $\phi^\alt_\bullet$ to
$\phi'^\alt_\bullet$.
\end{enumerate}
\end{lemma}
\begin{proof} Clearly it suffices to show (i). To this
aim, for every $i\in I$ let
$(\theta_i,\theta'_i):U_i\to U'_{\tau(i)}\times U'_{\tau'(i)}$
be the unique morphism whose composition with the projection
on the first (resp. second) factor equals $\theta_i$ (resp.
$\theta'_i$), and for every $n\in\N$, every $k=0,\dots,n$,
and every $\underline t\in I^{n+1}$, set
$$
\begin{aligned}
\tau(k,\underline t):=\,& (\tau(t_0),\dots,\tau(t_k),
\tau'(t_k),\dots,\tau'(t_n))\in I^{n+2} \\
\theta^{(k)}_{\underline t}:=\,& \theta_{t_0}\times\cdots
\times\theta_{t_{k-1}}\times(\theta_{t_k},\theta'_{t_k})
\times\theta'_{t_{k+1}}\times\cdots\times\theta'_{t_n}:
U_{\underline t}\to U'_{\tau(k,\underline t)}.
\end{aligned}
$$
Then $\theta^{(k)}_{\underline t}$ induces a natural morphism
$h^{(k)}_{\underline t}:j_{\underline t!}\Z_{U_{\underline t}}\to
j_{\tau(k,\underline t)!}\Z_{U'_{\tau(k,\underline t)}}$
and we set
$$
h_{\underline t}:=\sum_{k=0}^n (-1)^k\cdot h^{(k)}_{\underline t}:
j_{\underline t!}\Z_{U_{\underline t}}\to R_{n+1}(\fU').
$$
Clearly there exists a unique morphism of $\Z_T$-modules
$h_n:R_n(\fU)\to R_{n+1}(\fU')$ whose restriction to
$j_{\underline t!}\Z_{U_{\underline t}}$ agrees with $h_{\underline t}$
for every $\underline t\in I^{n+1}$. Let also
$h_{-1}:R_{-1}(\fU)\to R_0(\fU')$ be the zero morphism; a
direct calculation shows that the system $(h_{n-1}~|~n\in\N)$
yields the sought homotopy : the details shall be left
to the reader.
\end{proof}

\begin{remark}\label{rem_explicit-homotopy}
(i)\ \ 
In the situation of \eqref{subsec_alt-refine-map}, consider
the case where $U_{i_0}=1_T$ for some $i_0\in I$. Then $\fU$ can
be refined by itself in two different ways : we may take for
$\tau:I\to I$ the identity map (with $\theta_i:=\one_{U_i}$
for every $i\in I$), and we have also the map $\tau':I\to I$
such that $\tau'(i)=i_0$ for every $i\in I$ (and then there
is a unique choice of $\theta_i$). Obviously, the morphism
$\phi_\bullet:R^\alt_\bullet(\fU)\to R^\alt_\bullet(\fU)$ associated
with $\tau$ is the identity map. For the morphism
$\phi'_\bullet:R^\alt_\bullet(\fU)\to R^\alt_\bullet(\fU)$ associated
with $\tau'$ it is easily seen that
$$
\phi'_{-1}=\one_{\Z_T}
\qquad\text{and}\qquad
\phi'_n=0
\qquad
\text{for every $n>0$}
$$
and $\phi'_0$ is the sum of the natural morphisms
$j_{t!}\Z_{U_t}\to j_{i_0!}\Z_{U_{i_0}}=\Z_T$. The homotopy
$h_\bullet$ from $\phi_\bullet$ to $\phi'_\bullet$ furnished
by lemma \ref{lem_mad-combinatorics} can be described as
follows. First, we may suppose that $i_0$ is the maximal
element of $I$. A simple inspection shows that
$h_{-1}:R_{-1}^\alt(\fU)\to R_0^\alt(\fU)$ is the zero map.
Next, for every $n\in\N$ and every
$\underline t\in I^{n+1}_\alt$, notice that
$U_{\underline t}=U_{(\underline t,i_0)}$, where $(\underline t,i_0)$
denotes the sequence $(t_0,\dots,t_n,i_0)$; it follows that
we may write
$$
h_n=\sum_{\underline t\in I^{n+1}_\alt}\rho(\underline t)\cdot
h_{\underline t}
$$
where $h_{\underline t}:j_{\underline t}\Z_{U_{\underline t}}\to
j_{(\underline t,i_0)!}\Z_{U_{(\underline t,i_0)}}$ is the identity map,
and $\rho(\underline t)$ equals $0$ if $t_n=i_0$, and
equals $1$ otherwise.

(ii)\ \
In the situation of (i), we deduce a homotopy from
$\one_{R^\alt_\bullet(\fU)}$ to the zero map. Indeed,
let $h'_{-1}:\Z_T\to R_0(\fU)$ be the composition of
the identity map $\Z_T\to j_{i_0!}\Z_{U_{i_0}}$ and the
inclusion map $j_{i_0!}\Z_{U_{i_0}}\to R_0(\fU)$. Let
also $h'_n:=h_n$ for every $n\in\N$. A simple
inspection shows that
$$
\phi'_0=h'_{-1}\circ d_0
$$
(where $d_0:R_0(\fU)\to\Z_T$ is the differential of
$R_\bullet(\fU)$). It follows easily that
$(h'_{n-1}~|~n\in\N)$ is the sought homotopy.

(iii)\ \
Let $\cF$ be any $\Z_T$-module. The homotopy $h'_\bullet$
of (i) induces a homotopy from the identity automorphism
of the augmented \v{C}ech complex $C^\bullet_\alt(\fU,\cF)$
to its zero endomorphism. Under the natural identifications
of remark \ref{rem_justify-name}(ii), this homotopy $h^\bullet$
can be described as follows. 
\begin{itemize}
\item
$h^0:C^0(\fU,\cF)\to\Gamma(\cF)$ is the projection
onto the factor $\cF(U_{i_0})=\Gamma(\cF)$.
\item
For every $n\geq 0$, the map $h^{n+1}:C^{n+1}(\fU,\cF)\to C^n(\fU,\cF)$
is the projection that sends to zero every factor
$\cF(U_{t_0,\dots,t_n})$ with $t_n<i_0$, and whose restriction
to every factor $\cF(U_{t_0,\dots,t_{n-1},i_0})$ is the identity
map $\cF(U_{t_0,\dots,t_{n-1},i_0})\isom\cF(U_{t_0,\dots,t_{n-1}})$.
\end{itemize}
\end{remark}

\begin{remark}\label{rem_double-cechs}
(i)\ \
In the situation of \eqref{subsec_choose-total}, it
is clear that the rule $\cF\mapsto R^\alt_\bullet(\fU,\cF)$
defines a functor from $\Z_T$-modules to chain complexes
of $\Z_T$-modules. The latter then extends naturally to
a functor from the category of chain complexes of
$\Z_T$-modules to the category of double complexes
of $\Z_T$-modules :
$$
\sC(\Z_T\Mod)\to\sC_2(\Z_T\Mod)
\qquad
\cF_\bullet\mapsto R^\alt_\bullet(\fU,\cF_\bullet).
$$
Likewise, $\perp_\bullet^\fU$ extends to a functor from
augmented simplicial $\Z_T$-modules to augmented
bisimplicial $\Z_T$-modules :
$$
\widehat{s}.\Z_T\Mod\to\widehat{s}.(\widehat{s}.\Z_T\Mod)
\qquad
\cF_\bullet\mapsto\perp^\fU_\bullet\cF_\bullet.
$$

(ii)\ \
Let now $\fU':=(U'_{i'}~|~i'\in I')$ be another family of
objects of $T$. Then we may consider the augmented
bisimplicial $\Z_T$-module
$$
\perp_{\bullet\bullet}^{\fU,\fU'}:=\perp^{\fU'}_\bullet(\perp^\fU_\bullet\Z_T)
$$
and we notice that
$$
(\perp_{\bullet\bullet}^{\fU,\fU'})^\Delta=\perp_\bullet^{\fU\times\fU'}
\qquad
\text{where
$\fU\times\fU':=(U_i\times U'_{i'}~|~(i,i')\in I\times I')$}.
$$
In case both $U_i$ and $U'_{i'}$ are subobjects of $1_T$ for
every $i\in I$ and $i'\in I'$, we may also define the chain
double complex
$$
R^\alt_{\bullet\bullet}(\fU,\fU'):=R^\alt_\bullet(\fU',(R^\alt_\bullet(\fU))
$$
and we notice that (notation of example \ref{ex_monoidal}(i))
$$
R^\alt_{\bullet\bullet}(\fU,\fU')=
R^\alt_\bullet(\fU)\boxtimes_A R^\alt_\bullet(\fU').
$$
\end{remark}

\sset\subsubsection{}\label{subsec_pseudo-Leray}
In the situation of \eqref{subsec_Cech-resolution}, let
$\cA$ be any $T$-ring; we consider the functor
$$
\cA\Mod\to\sC(\Gamma(\cA)\Mod)
\qquad
\cF\mapsto\bar C{}^\bullet(\fU,\cF):=
\Hom_\Z(\bar R_\bullet(\fU),\cF[0])
$$
that assigns to every $\cA$-module $\cF$ its
{\em \v{C}ech complex}, and we define as well the
{\em \v{C}ech cohomology functor of $\cF$} relative
to the covering $\fU$, by setting :
$$
H^i(\fU,\cF):=H^i\bar C{}^\bullet(\fU,\cF)
\qquad
\text{for every $i\in\N$}.
$$
Suppose now that $\fU$ is a covering of $T$; since $\cF$
represents a sheaf for the canonical topology on $T$, we
easily see that the system of restriction maps
$(\Gamma(\cF)\to\cF(U_i)~|~i\in I)$ induces a natural
isomorphism
\set\begin{equation}\label{eq_goes-opposite-way}
\Gamma(\cF)\isom H^0(\fU,\cF).
\end{equation}
On the other hand, let $K^\bullet$ be any bounded below
complex of $\cA$-modules; lemma
\ref{lem_cotriples-are-back}(iii) yields natural isomorphisms
in $\sD^+(\cA\Mod)$ and in $\sD^+(\Gamma(\cA)\Mod)$
$$
R\cHom^\bullet_\Z(\bar R{}_\bullet(\fU),K^\bullet)\isom
R\underline\Gamma K^\bullet
\qquad
R\Hom^\bullet_\Z(\bar R{}_\bullet(\fU),K^\bullet)\isom
R\Gamma K^\bullet.
$$
Hence, after fixing an injective resolution $K^\bullet\isom\cI^\bullet$
we get a natural isomorphism
$$
R\Gamma K^\bullet\isom
\Tot\,\bar C{}^\bullet(\fU,\cI^\bullet)
\qquad
\text{in $\sD^+(\Gamma(\cA)\Mod)$}.
$$
Thus, for every $p\in\Z$ let also define a presheaf
$\cH^p(K^\bullet)$ on $T$ by the rule
$$
\cH^p(K^\bullet)(U):=R^p\Gamma(U,K^\bullet)
\qquad
\text{for every object $U$ of $T$}.
$$
We may regard $\cH^p(K^\bullet)$ as an object of the $\sV$-topos
$T^\wedge_\sV$, for some universe $\sV$ such that $T$ is $\sV$-small,
and the family $\fU$ can be regarded as a system of objects
of $T^\wedge_\sV$, via the Yoneda imbedding $T\to T^\wedge_\sV$.
Then, the \v{C}ech cohomology functors of $\cH^p(K^\bullet)$
relative to the covering $\fU$ are also well defined, and
we deduce a $2$-spectral sequence
\set\begin{equation}\label{eq_spectral-Leray}
E(\fU)_2^{p,q}:=H^p(\fU,\cH^q(K^\bullet))\Rightarrow
R^{p+q}\Gamma K^\bullet.
\end{equation}
Notice also that
$$
\cH^0(\cF[0])=\cF
\qquad
\text{for every $\Z_T$-module $\cF$}
$$
whence, natural maps
$$
\Psi^p(\fU,\cF):H^p(\fU,\cF)\to R^p\Gamma\cF
\qquad
\text{for every $p\in\N$}
$$
such that $\Psi^0(\fU,\cF)$ is always an isomorphism
(in fact, it is the inverse of the isomorphism
\eqref{eq_goes-opposite-way}), and $\Psi^1(\fU,\cF)$
is always an injective map.

\sset\subsubsection{}\label{subsec_altern-pseudo-Leray}
In case $U_i$ is a subobject of $1_T$ for every $i\in I$,
we may repeat the considerations of \eqref{subsec_pseudo-Leray}
with the {\em alternating \v{C}ech complex}
$$
\bar C{}^\bullet_\mathrm{alt}(\fU,\cF):=
\Hom_\Z(\bar R{}^\alt_\bullet(\fU),\cF[0]).
$$
Namely, let us set
$$
H^i_\alt(\fU,\cF):=H^i\bar C{}_\alt^\bullet(\fU,\cF)
\qquad
\text{for every $i\in\N$}.
$$
Then, proposition \ref{prop_Cech-resolve}(ii) (together
with example \ref{ex_hot-morphisms}(ii)) yields a natural
isomorphism
$$
H^i_\alt(\fU,\cF)\isom H^i(\fU,\cF)
\qquad
\text{for every $i\in\N$}
$$
and if $\fU$ is a covering of $T$, obviously we get
therefore a natural isomorphism
\set\begin{equation}\label{eq_same-but-altern}
R\Gamma K^\bullet\isom
\Tot\,\bar C{}_\alt^\bullet(\fU,\cI^\bullet)
\qquad
\text{in $\sD^+(\Gamma(\cA)\Mod)$}.
\end{equation}
for any injective resolution $K^\bullet\isom\cI^\bullet$, and
a $2$-spectral sequence
\set\begin{equation}\label{eq_alternate-spectral-Leray}
F(\fU)_2^{p,q}:=H^p_\alt(\fU,\cH^q(K^\bullet))\Rightarrow
R^{p+q}\Gamma K^\bullet.
\end{equation}
Taking $K^\bullet:=\cF[0]$, there follows a natural map
$$
\Psi^p_\alt(\fU,\cF):H^p_\alt(\fU,\cF)\to R^p\Gamma\cF
\qquad
\text{for every $p\in\N$}
$$
and again, $\Psi^0_\alt(\fU,\cF)$ is an isomorphism,
whereas $\Psi^1_\alt(\fU,\cF)$ is always an injective
map. The following immediate corollary is often useful :

\begin{corollary}\label{cor_Leray}
In the situation of \eqref{subsec_pseudo-Leray},
suppose furthermore that
$$
\cH^q(\cF[0])(U_{\underline t})=0
\qquad
\text{for every integer $q>0$ and every $\underline t\in I^q$}.
$$
Then we have :
\begin{enumerate}
\item
The map $\Psi^p(\fU,\cF)$ is an isomorphism for every $p\in\N$.
\item
If $U_i$ is a subobject of $1_T$ for every $i\in\N$, then
the map $\Psi^p_\alt(\fU,\cF)$ is an isomorphism for every
$p\in\N$.
\end{enumerate}
\end{corollary}
\begin{proof} The assumption of (i) (resp. (ii)) implies that
$E(\fU)^{p,q}_2=0$ (resp. $F(\fU)^{p,q}_2=0$) for every $q>0$,
whence both contentions.
\end{proof}

\sset\subsubsection{}\label{subsec_go-to-lim-and-Cech}
Let now $C:=(\cC,J)$ be a small site whose finite non-empty
products and fibre products are representable. Then, the
finite products of $\cC/X$ are representable, for every
$X\in\Ob(\cC)$; for every such $X$, we denote by $J_X$
the topology on $\cC/X$ induced by $J$ (see
\eqref{subsec_topol-on-C-over-X}).

Now, let $S$ be either $C$ or $(\cC/X,J_X)$ for some
$X\in\Ob(\cC)$, and denote by $\cJ(S)$ the set of all
families of objects of $S$ that generate a sieve covering
the final object of $T:=S^\sim$, for the canonical topology
on $T$. Then $\cJ(S)$ is endowed with a partial ordering,
by declaring that $\fU\leq\fU'$ if and only if the sieve
generated by $\fU$ contains the sieve generated by $\fU'$,
for every $\fU,\fU'\in\cJ(S)$.
Notice that this is equivalent to saying that the family
$\fU'$ is a refinement of $\fU$, when we regard $\fU$ and
$\fU'$ as families of objects of the topos $T$.
Especially, when $\fU\leq\fU'$, the discussion of
\eqref{subsec_maddening-combinatorics} yields a map of
complexes $\phi_{\fU,\fU'}:R_\bullet(\fU')\to R_\bullet(\fU)$,
and lemma \ref{lem_mad-combinatorics} says that the
homotopy class of $\phi_{\fU,\fU'}$ depends only on $\fU$
and $\fU'$, so we get a well defined functor
$$
R_\bullet:\cJ(S)^o\to\Hot(\Z_T\Mod)
\qquad
\fU\mapsto R_\bullet(\fU).
$$
Now, let $\cF$ be any abelian sheaf on $S$. We
then get, for every $i\in\N$, an induced functor
$$
H^i(-,\cF):\cJ(S)\to\Z\Mod
\qquad
\fU\mapsto H^i(\fU,\cF).
$$
Notice now that $\cJ(S)$ is a small and filtered
category; we may then define
$$
\check{H}^i(S,\cF):=\colim_{\fU\in\cJ(S)}H^i(\fU,\cF)
$$
and clearly the rule $\cF\mapsto\check{H}^i(S,\cF)$ yields
a well defined functor
$$
\check{H}^i(S,-):\Z_T\Mod\to\Z\Mod
\qquad
\text{for every $i\in\N$}.
$$
In case $S=(\cC/X,J_X)$, we shall denote this functor by
$\check{H}(X,-)$. By the same token, it is also clear that
we get as well a natural morphism of $2$-spectral sequences
$$
E(\fU)_2^{\bullet\bullet}\to E(\fU')_2^{\bullet\bullet}
\qquad
\text{whenever $\fU\leq\fU'$}
$$
for every bounded below complex $K^\bullet$ of $\Z_T$-modules
whence, after taking colimits, a $2$-spectral sequence
$$
E_2^{p,q}:=\check{H}^p(S,\cH^q(K^\bullet))\Rightarrow
R^{p+q}\Gamma K^\bullet
$$
that, in turn, gives us again natural maps, for every
$\Z_T$-module $\cF$
$$
\check{\Psi}{}^p(S,\cF):\check{H}^p(S,\cF)\to R^p\Gamma\cF
\qquad
\text{for every $p\in\N$}
$$
which, in case $S=(\cC/X,J_X)$, shall be denoted simply
$\check{\Psi}{}^p(X,\cF)$.

\begin{lemma}\label{lem_degenerate-Cechist}
In the situation of \eqref{subsec_go-to-lim-and-Cech},
let $\cF$ be any $\Z_T$-module, and take $K^\bullet:=\cF[0]$.
Then we have :
\begin{enumerate}
\item
$E^{0,q}_2=0$ for every $q>0$.
\item
Both $\check{\Psi}{}^0(S,\cF)$ and $\check{\Psi}{}^1(S,\cF)$
are isomorphisms, and $\check{\Psi}{}^2(S,\cF)$ is injective.
\end{enumerate}
\end{lemma}
\begin{proof}(i): Fix any injective resolution $\cF\to\cI^\bullet$,
and let $\bar s\in E^{0,q}_2$ be any element; by definition,
there exists a family $(U_i~|~i\in I)$ of objects of $S$ that
generates a sieve covering the final object of $S^\sim$, such
that $\bar s$ is represented by an element
$$
\bar s_\bullet\in\Ker\Bigl(\prod_{i\in I}H^q(\cI^\bullet(U_i))
\xrightarrow{\ d_v^{0,q}\ }
\prod_{(i,j)\in I^2}H^q(\cI^\bullet(U_{(i,j)})\Bigr)
$$
where $d_v^{0,q}$ is the differential of the \v{C}ech complex in
degree $0$, as in remark \ref{rem_justify-name}(i). In turn,
the class $\bar s_i\in H^q(\cI^\bullet(U_i))$ is represented
by a section $s_i\in\cI^q(U_i)$ for every $i\in I$, such that
$d^{0,q}_h(s_i)=0$, where $d^{0,q}_h:\cI^q(U_i)\to\cI^{q+1}(U_i)$
denotes the differential in degree $q$ of the complex $\cI^\bullet$.
Since the latter is exact in every degree $q>0$, it follows
that for every $i\in I$ there exist a family
$(U_{i,\lambda}\to U_i~|~\lambda\in\Lambda_i)$ of objects of
$\cC/U_i$ that generates a sieve covering the final object of
$(\cC/U_i,J_{U_i})$, and for every $\lambda\in\Lambda_i$ a
section $s_{i,\lambda}\in\cI^{q-1}(U_{i,\lambda})$ such that
$d^{0,q-1}_h(s_{i,\lambda})=s_{i|U_{i,\lambda}}$. Set
$\fU:=\bigcup_{i\in I}\{U_{i,\lambda}~|~\lambda\in\Lambda_i\}$;
then $\fU$ lies in $\cJ(S)$, and it is easily seen that the
image of $\bar s$ vanishes in $H^0(\fU,\cH^q(\cF))$, whence
the contention.

(ii) is an immediate consequence of (i).
\end{proof}

\begin{theorem}[Cartan]\label{th_Cartan}
In the situation of \eqref{subsec_go-to-lim-and-Cech}, let
$\cF$ be any $\Z_T$-module such that
$$
\check{H}{}^i(X,\cF)=0
\qquad
\text{for every $X\in\Ob(\cC)$ and every $i>0$}.
$$
Then we have :
\begin{enumerate}
\item
The map $\check{\Psi}^q(C,\cF)$ is an isomorphism
for every $q\in\N$.
\item
More precisely, the map $\Psi^q(\fU,\cF)$ is an isomorphism
for every $q\in\N$ and every family $\fU$ of objects of\/
$\cC$ that covers the final object of\/ $T$.
\end{enumerate}
\end{theorem}
\begin{proof}(i): We argue by induction on $q\in\N$, and
notice that the assertion for $q\leq 1$ is already known
without any assumption on $\cF$, by lemma
\ref{lem_degenerate-Cechist}(ii). Thus, let $i\geq 2$ and
suppose that the assertion of the theorem is already
known for every $q<i$, every site $C$ and every abelian
sheaf $\cF$ on $C$. Notice that if $X\in\Ob(\cC)$ and
$\phi:Y\to X$ is any object of $\cC/X$, we have a natural
isomorphism of categories
$$
(\cC/X)/\phi\isom\cC/Y
$$
and the topology $J_\phi$ induced on $(\cC/X)/\phi$ by $J_X$
agrees, under this identification, with the topology $J_Y$.
Therefore, the site $\cC/X$ fulfills the assumptions of
the theorem, and by inductive assumption we deduce that
$\check{\Psi}^q(X,\cF)$ is an isomorphism for every $q<i$,
and therefore $R^q\Gamma(X,\cF)=0$ for every $X\in\Ob(\cC)$
and every $q=1,\dots,i-1$. Summing up, we get $\cH^q(\cF)=0$
whenever $1\leq q<i$, whence $E^{pq}_2=0$ for every
$q=1,\dots,i-1$ (notation of \eqref{subsec_go-to-lim-and-Cech}).
We also know that $E^{0q}_2=0$ for every $q>0$, by lemma
\ref{lem_degenerate-Cechist}(i); in this situation, it is
easily seen that $E_2^{i,0}=E^{i0}_\infty$ and $E^{i-p,p}_\infty=0$
for every $p>0$, and then $\check{\Psi}^i(X,\cF)$ is indeed
an isomorphism.

(ii) is an immediate consequence of (i) and corollary
\ref{cor_Leray}(i).
\end{proof}

\sset\subsubsection{}\label{subsec_cohom-of-affine}
We wish now to apply the foregoing results to the cohomology
of quasi-coherent modules on a scheme. To begin with, let
$A$ be a ring, $M$ an $A$-module; set $X:=\Spec\,A$ and
denote by $\cM$ the quasi-coherent $\cO_{\!X}$-module arising
from $M$. Let also $\bff:=(f_1,\dots,f_r)$ be a sequence of
elements of $A$, and for every integer $n\geq -1$ and every
injective order-preserving map
$$
\underline t:[n]\to\Sigma:=\{1,\dots,r\}
\qquad
k\mapsto t_k
$$
let $A_{\underline t}:=A[f_{t_0}^{-1}\cdots f_{t_n}^{-1}]$ and
$U_{\underline t}:=\Spec\,A_{\underline t}$; we can describe as
follows the $A_{\underline t}$-module $\cM(U_{\underline t})$.
For every such $\underline t$, consider the system of
$A$-modules $((M^{(k)}_{\underline t},\phi_k)~|~k\in\N)$
such that $M^{(k)}_{\underline t}:=M$ for every $k\in\N$, and 
the map $\phi_k:M^{(k)}_{\underline t}\to M^{(k+1)}_{\underline t}$
is the scalar multiplication by $f_{t_0}\cdots f_{t_n}$ if
$n\geq 0$, and it is $\one_M$ when $n=-1$.
 Then we have a natural identification
\set\begin{equation}\label{eq_sections-as-colims}
\cM(U_{\underline t})\isom\colim_{k\in\N}\,M^{(k)}_{\underline t}
\qquad
\text{for every $n\geq-1$ and every
$\underline t:[n]\to\Sigma$}.
\end{equation}
Moreover, for every face map $\eps_i:[n-1]\to[n]$ in
$\Delta$, we can describe as follows the restriction map
$\rho_{i,\underline t}:\cM(U_{\underline t\circ\eps_i})\to
\cM(U_{\underline t})$. We consider the morphism of directed
systems
$$
\rho^{(\bullet)}_i:M^{(\bullet)}_{\eps_i^*(\underline t)}\to
M^{(\bullet)}_{\underline t}
\qquad
\text{such that}
\qquad
\rho^{(k)}_i:=f^k_{t_i}\cdot\one_M
\quad
\text{for every $k\in\N$}.
$$
Then, under the identifications \eqref{eq_sections-as-colims},
the map $\rho_{i,\underline t}$ corresponds to the map
$$
\colim_{k\in\N}\rho^{(k)}_i:\colim_{k\in\N}\,M^{(k)}_{\underline t\circ\eps_i)}
\to\colim_{k\in\N}\,M^{(k)}_{\underline t}.
$$
For every $n\geq-1$ let $\Sigma^{n+1}_\alt$ be the set of
all injective order-preserving maps $[n]\to\Sigma$, and set
$$
D^n_{(k)}:=\Hom_\Set(\Sigma^{n+1}_\alt,M)
\qquad
\text{for every $k\in\N$}
$$
which we endow with the $A$-module structure inherited from
$M$; consider also the map
$$
d^n_{(k)}:D^n_{(k)}\to D^{n+1}_{(k)}
\qquad
(\mu:\Sigma^{n+1}_\alt\to M)\mapsto
\Bigl(
\underline t\mapsto\sum_{i=0}^{n+1}(-1)^i\cdot f_{t_i}^k\cdot
\mu(\underline t\circ\eps_i)
\Bigr)
$$
for every $k\in\N$. Then the system $(D^\bullet_{(k)},d^\bullet_{(k)})$
is a well defined complex of $A$-modules, for every $k\in\N$.
Lastly, for every $k\in\N$ we consider the morphism of complexes
$$
D^\bullet_{(k)}\to D^\bullet_{(k+1)}
\qquad
(\mu:\Sigma^{n+1}_\alt\to M)\mapsto
(\underline t\mapsto f_{t_0}\cdots f_{t_n}\cdot\mu(\underline t)).
$$
Summing up, and comparing with remark \ref{rem_justify-name}(i)
we obtain a natural identification
$$
C^\bullet_\alt(\fU,\cM)\isom\colim_{k\in\N}D^\bullet_{(k)}
\qquad
\text{where $\fU:=(U_1,\dots,U_r)$}.
$$
On the other hand, notice that, for every $n\geq-1$,
the free $A$-module $\Lambda^{n+1}_AA^{\oplus r}$ admits
the standard basis $(e_{t_0}\wedge\cdots\wedge e_{t_n}~|~
\underline t\in\Sigma_\alt^{n+1})$, so we have as well a
natural identification
$$
\Hom_A(\Lambda^{n+1}_AA^{\oplus r},M)\isom D^n_{(k)}
\qquad
\text{for every $n\geq-1$}
$$
and a direct inspection reveals that, under this latter
identification, the differential $d^n_{(k)}$ of
$D^\bullet_{(k)}$ corresponds to the differential
$\Hom_A(d_{\bff^k,n},M)$, where $\bff^k:=(f_1^k,\dots,f^k_r)$,
and $d_{\bff^k,n}$ is the differential in degree $n$ of the
Koszul complex of the sequence $\bff^k$ (see remark
\ref{rem_koszul-alg}(ii)). Thus, we get finally a natural
isomorphism of complexes of $A$-modules
\set\begin{equation}\label{eq_Koszul-Cech}
C^\bullet_\alt(\fU,\cM)\isom\colim_{k\in\N}\bK^\bullet(\bff^k,M)[1]
\end{equation}
where the transition maps
$\bK^\bullet(\bff^k,M)\to\bK^\bullet(\bff^{k+1},M)$ are
the morphisms $\bphi_\bff$ of \eqref{subsec_def-bphis}.
As a corollary, we get the following classical result of
Grothendieck :

\begin{theorem}\label{th_Cech-resolve}
In the situation of \eqref{subsec_cohom-of-affine},
the following holds :
\begin{enumerate}
\item
The natural map $M[0]\to R\Gamma\,\cM$ is an isomorphism
in $\sD(\cO_{\!X}\Mod)$.
\item
Let $Y$ be any separated scheme, $\cF^\bullet$ a bounded below
complex of quasi-coherent $\cO_Y$-modules, and $\fU$ any
affine open covering of\/ $Y$. Then the natural map
$$
\Tot\,\bar C{}_\mathrm{alt}^\bullet(\fU,\cF^\bullet)\to
R\Gamma\cF^\bullet
$$
is an isomorphism in $\sD^+(\cO_Y(Y)\Mod)$.
\item
Especially, if\/ $Y$ is a separated scheme that can be covered
by $r$ affine open subsets, then $H^i(Y,\cF)=0$ for every
$i\geq r$, and every quasi-coherent $\cO_Y$-module $\cF$.
\end{enumerate}
\end{theorem}
\begin{proof}(i): Since the affine open subsets of $X$ of the
form $\Spec\,A[f^{-1}]$ for $f\in A$ are a basis of the
topology of $X$ that is closed under finite intersections,
theorem \ref{th_Cartan} reduces to showing that
$$
H^0(\fU,\cM)=M
\qquad\text{and}\qquad
H^i(\fU,\cM)=0
\qquad
\text{for every $i>0$}
$$
where $\fU:=(U_i,\dots,U_r)$ is the finite affine open covering
of $X$ associated with a sequence $\bff:=(f_1,\dots,f_r)$ as
in \eqref{subsec_cohom-of-affine}. However, the condition
$X=\bigcup_{i=1}^rU_i$ is equivalent to saying that the ideal
generated by the system $\bff$ is $A$, in which case the
same holds for the ideal generated by $\bff^k$, for every
$k\in\N$, and then lemma \ref{lem_koszul-vanish}(iii) says
that the complex $H^\bullet(\bff^k,M)$ is homotopically
trivial for every $k\in\N$. Taking into account
\eqref{eq_Koszul-Cech}, the contention follows.

(ii): Pick any resolution $\cF^\bullet\isom\cI^\bullet$ by a
bounded below complex of injective $\cO_Y$-modules; in view
of \eqref{eq_same-but-altern}, it suffices to show that the
natural map of double complexes
$$
\bar C{}_\mathrm{alt}^\bullet(\fU,\cF^\bullet)\to
\bar C{}_\mathrm{alt}^\bullet(\fU,\cI^\bullet)
$$
induces a quasi-isomorphism on total complexes. However, notice
that the open subset $U_{\underline t}\subset Y$ is affine, for
every $n\in\N$ and every $\underline t\in I^{n+1}_\alt$ (notation
of \eqref{subsec_choose-total}); hence, it suffices to check
that the induced map $\cF^\bullet(V)\to\cI^\bullet(V)$ is a
quasi-isomorphism, for every affine open subset $V\subset Y$.
To this aim, consider the spectral sequence
$$
E_1^{p,q}:=H^p(V,\cF^q)\Rightarrow H^{p+q}\cI^\bullet(V).
$$
Since $V$ is affine, we have $E_1^{p,q}=0$ for every $p>0$,
due to (i); on the other hand, the differential
$d_1^{0,q}:E_1^{0,q}\to E_1^{0,q+1}$ is nothing else than the
differential $d^q(V):\cF^q(V)\to\cF^{q+1}(V)$, for every
$q\in\Z$, whence the claim.

(iii): Indeed, if $\fU$ is such a covering, then the
complex $\bar C{}^\bullet_\alt(\fU,\cF)$ is an object
of the category $\sC^{[0,r-1]}(\cO_Y(Y)\Mod)$, so the
assertion follows from (ii).
\end{proof}

We conclude this section with an application of \v{C}ech
cohomology to the computation of sheaf cohomology after
change of base scheme; much more can be found in
\cite[Ch.III, \S6]{EGAIII-2}.

\sset\subsubsection{}\label{subsec_local-Tors}
Namely, consider a diagram of schemes
$$
\xymatrix{ X_2 \ar[r] \ar[d]_{\pi_2} &
X_1 \ar[r] \ar[d]_{\pi_1} & X_0 \ar[d]^{\pi_0} \\
S_2 \ar[r] & S_1 \ar[r]^-\phi & S_0
}$$
whose two square subdiagram are cartesian, and bounded
above complexes of quasi-coherent $\cO_{\!X_0}$-modules
$\cF_0^\bullet$ and quasi-coherent $\cO_{\!S_i}$-modules
$\cG^\bullet_i$ for $i=1,2$. Then, for every $q\in\N$
there exists a quasi-coherent $\cO_{\!X_1}$-module
$$
\cTor_q^{S_0}(\cG^\bullet_1,\cF^\bullet_0)
$$
characterized, up to unique isomorphism, by the following
properties (\cite[Ch.III, \S6.5]{EGAIII-2}) :
\begin{enumerate}
\alphaenu
\item
For every affine open subsets $V_0\subset S_0$,
$V\subset S_1\times_{S_0}V_0$ and $U\subset X_0\times_{S_0}V_0$,
there exists an isomorphism of $\cO_{\!X_1}(V\times_{S_0}U)$-modules
$$
\omega_{V,U}:
\Gamma(V\times_{S_0}U,\cTor_q^{S_0}(\cG^\bullet_1,\cF^\bullet_0))
\isom
H_q(\cG^\bullet_1(V)\derotimes_{\cO_{\!S_0}(V_0)}\cF^\bullet_0(U)).
$$
\item
For every $V_0$, $V$, $U$ as in (a) and every inclusion of
affine open subsets $V'_0\subset V_0$,
$V'\subset V\cap(S_1\times_{S_0}V'_0)$ and
$U'\subset U\cap(X_0\times_{S_0}V'_0)$, the resulting diagram
commutes :
$$
\xymatrix{
\Gamma(V\times_{S_0}U,\cTor_q^{S_0}(\cG_1^\bullet,\cF^\bullet_0))
\ar[rr]^-{\omega_{V,U}} \ar[d] & &
H_q(\cG^\bullet_1(V)\derotimes_{\cO_{\!S_0}(V_0)}\cF^\bullet_0(U)) \ar[d] \\
\Gamma(V'\times_{S_0}U',\cTor_q^{S_0}(\cG^\bullet_1,\cF^\bullet_0))
\ar[rr]^-{\omega_{V',U'}} & &
H_q(\cG^\bullet_1(V')\derotimes_{\cO_{\!S_0}(V'_0)}\cF^\bullet_0(U'))
}$$
whose left vertical arrow is the restriction map of the
sheaf $\cTor_q^{S_0}(\cG_1^\bullet,\cF^\bullet_0)$, and whose right
vertical arrow is the map $H_q(\rho'\derotimes_\rho\rho'')$
associated with the restriction maps
$\rho:\cO_{\!S_0}(V_0)\to\cO_{\!S_0}(V'_0)$,
$\rho':\cG^\bullet(V)\to\cG^\bullet(V')$ and
$\rho'':\cF^\bullet_0(U)\to\cF^\bullet_0(U')$.
\end{enumerate}

\begin{proposition}\label{prop_local-Tors}
In the situation of \eqref{subsec_local-Tors}, the following
holds :
\begin{enumerate}
\item
There exists a natural (homological) $2$-spectral sequence
$$
\cTor_p^{S_1}(\cG_2^\bullet,\cTor_q^{S_0}(\cO_{\!S_1}[0],\cF^\bullet_0))
\Rightarrow
\cTor^{S_0}_{p+q}(\cG_2^\bullet,\cF^\bullet_0).
$$
\item
Suppose that $\pi_0$ is quasi-compact and separated, $\cF_0$
is a bounded complex, and both $S_1$ and $S_0$ are affine.
Set $A_i:=\cO_{\!S_i}(S_i)$ for $i=0,1$. Then there exist two
natural spectral sequences :
$$
\begin{aligned}
E_2^{pq}:=\,& H^{-p}(X_1,\cTor_q^{S_0}(\cG^\bullet,\cF_0^\bullet))
\Rightarrow
H_{p+q}(R\Gamma\cG^\bullet\derotimes_{A_0}R\Gamma\cF^\bullet_0) \\
F_2^{pq}:=\,&
H_p(R\Gamma\cG^\bullet\derotimes_{A_0}H^{-q}(X_0,\cF^\bullet_0)[0])
\Rightarrow
H_{p+q}(R\Gamma\cG^\bullet\derotimes_{A_0}R\Gamma\cF^\bullet_0).
\end{aligned}
$$
\end{enumerate}
\end{proposition}
\begin{proof}(i): This is obtained from the standard
change of base ring spectral sequence, which is natural
in all arguments, and therefore globalizes immediately
to the sheaf-theoretic situation considered here : the
details shall be left to the reader.

(ii): Under our assumptions, $X_0$ is quasi-compact
and separated, hence we may find a finite affine covering
$\fU$ of $X_0$, and by theorem \ref{th_Cech-resolve}(ii),
the alternating \v{C}ech complex for $\fU$ computes the
cohomology of $\cF^\bullet_\bullet$; especially the complex
$R\Gamma\cF_0^\bullet$ is bounded. Moreover, we may find a
bounded above Cartan-Eilenberg resolution
$\cL^\bullet\isom\phi_*\cG^\bullet$ consisting of quasi-coherent
flat $\cO_{\!S_0}$-modules, and
$H_i(R\Gamma\cG^\bullet\derotimes_{A_0}R\Gamma\cF^\bullet_0)$ is
naturally isomorphic to the quasi-coherent $\cO_{\!S_1}$-module
associated with the $A_1$-module
$$
H_i((\Tot\,\bar C{}^\bullet_\alt(\fU,\cF^\bullet_0))
\otimes_{A_0}\cL^\bullet(S_0))
$$
for every $i\in\Z$. Then the first spectral sequence
is the standard spectral sequence attached to the
double complex $(\Tot\,\bar C{}^\bullet_\alt(\fU,\cF^\bullet_0))
\boxtimes_{A_0}\cL^\bullet(S_0)$; similarly one obtains the
second spectral sequence : details left to the reader.
\end{proof}

\subsection{Quasi-coherent modules}\label{sec_various-O-mod}
For any scheme $X$, we denote by
$$
\cO_{\!X}\Mod
\qquad
\cO_{\!X}\Mod_\qcoh
\qquad
\cO_{\!X}\Mod_\coh
\qquad
\cO_{\!X}\Mod_\mathrm{lfft}
$$
the category of all (resp. of quasi-coherent, resp. of coherent,
resp. of locally free of finite type) $\cO_{\!X}$-modules.
Recall that there is a natural functor
\set\begin{equation}\label{eq_arising-qcoh}
\cO_{\!X}(X)\Mod\to\cO_{\!X}\Mod_\qcoh
\qquad
M\mapsto M^\sim
\end{equation}
that assigns to every $\cO_{\!X}(X)$-module the quasi-coherent
module $M^\sim$ such that
$$
M^\sim(U):=M\otimes_{\cO_{\!X}(X)}\cO_{\!X}(U)
\qquad
\text{for every affine open subset $U\subset X$}.
$$
This functor is an equivalence, if $X$ is affine. We notice :

\begin{corollary}\label{cor_Ext-loc=glob}
Let $A$ be a ring, $M^\bullet$ (resp. $N_\bullet$) a bounded
below (resp. above) complex of $A$-modules. Set $X:=\Spec\,A$
and denote by $M^{\bullet\sim}$ (resp. $N_\bullet^\sim$) the
associated complex of quasi-coherent $\cO_{\!X}$-modules.
We have :
\begin{enumerate}
\item
The natural map
$$
R\Hom^\bullet_A(N_\bullet,M^\bullet)\to
R\Hom^\bullet_{\cO_{\!X}}(N^\sim_\bullet,M^{\bullet\sim})
$$
is an isomorphism in $\sD^+(A\Mod)$.
\item
If $A$ is coherent and $N_\bullet$ is a complex of finitely
presented $A$-modules, then the natural morphism
$$
R\Hom^\bullet_A(N_\bullet,M^\bullet)^\sim\to
R\cHom_{\cO_{\!X}}^\bullet(N^\sim_\bullet,M^{\bullet\sim})
$$
is an isomorphism in $\sD(\cO_{\!X}\Mod)$.
\end{enumerate}
\end{corollary}
\begin{proof}(i): We apply the trivial duality theorem
\ref{th_trivial-dual} to the unique morphism
$$
f:(X,\cO_{\!X})\to(\{\mathrm{pt}\},A)
$$
of ringed spaces, where $\{\mathrm{pt}\}$ denotes the one-point
space. Since $f$ is flat and all quasi-coherent $\cO_{\!X}$-modules
are $f_*$-acyclic, the assertion follows easily.

(ii): Due to (i), it suffices to show the following
assertion. For every $f\in A$, the natural map
$$
R\Hom^\bullet_A(N_\bullet,M^\bullet)\otimes_AA_f\to
R\Hom^\bullet_{A_f}(N_{\bullet,f},M^\bullet_f)
$$
is an isomorphism in $\sD(A_f\Mod)$. To this aim,
notice that, since $A$ is coherent, we may find a
resolution $P_\bullet$ of $N_\bullet$ consisting of
free $A$-modules of finite type (details left to the
reader); then we come down to checking that the
natural map
$\Hom_A(P_i,M^j)\otimes_AA_f\to\Hom_{A_f}(P_{i,f},M^j_f)$
is an isomorphism, which is obvious.
\end{proof}

We denote by
$$
\sD(\cO_{\!X}\Mod)_\qcoh
\qquad
\sD(\cO_{\!X}\Mod)_\coh
$$
the full triangulated subcategory of $\sD(\cO_{\!X}\Mod)$
consisting of the complexes $K^\bullet$ such that
$H^iK^\bullet$ is a quasi-coherent (resp. coherent)
$\cO_{\!X}$-module for every $i\in\Z$. As usual, we shall
use also the variants $\sD^+(\cO_{\!X}\Mod)_\qcoh$ (resp.
$\sD^-(\cO_{\!X}\Mod)_\qcoh$, resp. $\sD^b(\cO_{\!X}\Mod)_\qcoh$,
resp. $\sD^{[a,b]}(\cO_{\!X}\Mod)_\qcoh$) consisting of all
objects of $\sD(\cO_{\!X}\Mod)_\qcoh$ whose cohomology vanishes
in sufficiently large negative degree (resp. sufficiently large
positive degree, resp. outside a bounded interval, resp. outside
the interval $[a,b]$), and likewise for the corresponding
subcategories of $\sD(\cO_{\!X}\Mod)_\coh$. Obviously,
\eqref{eq_arising-qcoh} induces a natural functor
$$
\sD(\cO_{\!X}(X)\Mod)\to\sD(\cO_{\!X}\Mod)_\qcoh
\qquad
M^\bullet\mapsto M^{\bullet\sim}.
$$

\begin{proposition}\label{prop_replace-prop.12.3.5}
Let $f:Y\to X$ be a flat morphism of schemes, with $X$
locally coherent (see definition {\em \ref{def_coh-schemes}(i)}),
and $K_\bullet\in\Ob(\sD^-(\cO_{\!X}\Mod)_\coh)$ any complex.
We have :
\begin{enumerate}
\item
$R\cHom^\bullet_{\cO_{\!X}}(K_\bullet,L^\bullet)\in
\Ob(\sD^+(\cO_{\!X}\Mod)_\coh)$ for every
$L^\bullet\in\Ob(\sD^+(\cO_{\!X}\Mod)_\coh)$.
\item
For every $L^\bullet\in\Ob(\sD^+(\cO_{\!X}\Mod))$
the natural morphism
\set\begin{equation}\label{eq_filter-this}
f^*R\cHom^\bullet_{\cO_{\!X}}(K_\bullet,L^\bullet)\to
R\cHom^\bullet_{\cO_{\!Y}}(f^*K_\bullet,f^*L^\bullet)
\end{equation}
is an isomorphism in $\sD^+(\cO_{\!X}\Mod))$.
\end{enumerate}
\end{proposition}
\begin{proof}(i): The assertion is local on $X$, so we may
assume that $X=\Spec\,A$, for some coherent ring $A$.
In this case, we may find a bounded below (resp. above)
complex $M^\bullet$ (resp. $N_\bullet$) of finitely
presented $A$-modules, with isomorphisms
$K^\bullet\isom M^{\bullet\sim}$, $L_\bullet\isom N_\bullet^\sim$
in $\sD(\cO_{\!X}\Mod)$. By virtue of corollary
\ref{cor_Ext-loc=glob}(ii), we are then reduced to
checking that $R\Hom^\bullet_A(N_\bullet,M^\bullet)$
lies in $\sD(A\Mod_\coh)$. But this is easily seen,
since we may find a resolution of $N_\bullet$ by a
bounded above complex of free $A$-modules of finite
type (details left to the reader).

(ii): Again, the question is local on $X$, so we may assume
that $K_\bullet$ admits a resolution by a bounded above
complex $P_\bullet$ of free $\cO_{\!X}$-modules of finite
type. Since $f$ is flat, we have natural convergent spectral
sequences
$$
\begin{aligned}
E^{pq}_1:=f^*R^p\cHom^\bullet_{\cO_{\!X}}(P_q,L^\bullet)
\Rightarrow\, & f^*R^{p+q}\cHom^\bullet_{\cO_{\!X}}(K_\bullet,L^\bullet) \\
F^{pq}_1:=R^p\cHom^\bullet_{\cO_{\!Y}}(f^*P_q,f^*L^\bullet)
\Rightarrow\, & R^{p+q}\cHom^\bullet_{\cO_{\!Y}}(f^*K_\bullet,f^*L^\bullet)
\end{aligned}
$$
as well as a morphism of spectral sequences
$E^{\bullet\bullet}_\bullet\to F^{\bullet\bullet}_\bullet$,
such that $H^\bullet\eqref{eq_filter-this}$ is a morphism
of filtered $\cO_Y$-modules, for the two finite filtrations
induced by these spectral sequences on their abutments.
It then suffices to check that the induced morphism
$E^{pq}_\infty\to F^{pq}_\infty$ is an isomorphism for
every $p,q\in\Z$, and this in turn will follow, if we
show that the morphism $E^{pq}_1\to F^{pq}_1$ is an
isomorphism; but the latter assertion is obvious.
\end{proof}

\begin{corollary}\label{cor_dir-im-and-colim}
Let $f:Y\to X$ be a quasi-compact and quasi-separated morphism
of schemes, $\cF$ a flat quasi-coherent $\cO_{\!X}$-module, and
$\cG$ any $\cO_Y$-module. Then the natural map
\set\begin{equation}\label{eq_dir-im-and-colim}
\cF\otimes_{\cO_{\!X}}Rf_*\cG\to Rf_*(f^*\cF\otimes_{\cO_Y}\cG)
\end{equation}
is an isomorphism in $\sD(\cO_{\!X}\Mod)$.
\end{corollary}
\begin{proof} The assertion is local on $X$, hence we may assume
that $X=\Spec\,A$ for some ring $A$, and $\cF=M^\sim$ for some
flat $A$-module $M$. By \cite[Ch.I, Th.1.2]{La}, $M$ is the colimit
of a filtered family of free $A$-modules of finite rank; in view of
proposition \ref{prop_dir-im-and-colim}(ii), we may then assume that
$M=A^{\oplus n}$ for some $n\geq 0$, in which case the assertion is
obvious.
\end{proof}

\begin{remark}\label{rem_dir-im-and-colim}
Notice that, for an affine morphism $f:Y\to X$, the
map \eqref{eq_dir-im-and-colim} is an isomorphism for
any quasi-coherent $\cO_{\!X}$-module $\cF$ and any
quasi-coherent $\cO_Y$-module $\cG$. The details shall
be left to the reader.
\end{remark}

\begin{corollary}\label{cor_base-change-where}
Consider a cartesian diagram of schemes
$$
\xymatrix{ Y' \ar[r]^-{g'} \ar[d]_{f'} & Y \ar[d]^f \\
           X' \ar[r]^-g & X
}$$
such that $f$ is quasi-compact and quasi-separated, and
$g$ is flat. Then the natural map
\set\begin{equation}\label{eq_where-is-this}
g^*Rf_*\cG\to Rf'_*g'{}^*\cG
\end{equation}
is an isomorphism in $\sD(\cO_{\!X'}\Mod)$, for every
$\cO_Y$-module $\cG$.
\end{corollary}
\begin{proof} We easily reduce to the case where both
$X$ and $X'$ are affine, hence $g$ is an affine morphism.
In this case, it suffices to show that
$g_*\eqref{eq_where-is-this}$ is an isomorphism.
However, we have an essentially commutative diagram
$$
\xymatrix{
g_*(g^*Rf_*\cG) \ar[rr]^-\alpha \ar[d] & &
Rf_*\cG\otimes_{\cO_{\!X}}g_*\cO_{\!X'} \ar[d] \\
g_*(Rf'_*g'{}^*\cG) \ar[r]^-\sim &
Rf_*g'_*(g'{}^*\cG) \ar[r]^-\beta &
Rf_*(\cG\otimes_{\cO_{\!Y}}g'_*\cO_{\!Y'})
}$$
whose left vertical arrow is $g_*\eqref{eq_where-is-this}$
and whose right vertical arrow is the natural isomorphism
provided by corollary \ref{cor_dir-im-and-colim} (applied
to the flat $\cO_{\!Y}$-module $\cF:=g_*\cO_{\!Y'}$). Also,
$\alpha$ and $\beta$ are the natural
maps obtained as in \cite[Ch.0, (5.4.10)]{EGAI}; it is
easily seen that these are isomorphisms, for any affine
morphism $g$ and $g'$. The assertion follows.
\end{proof}

\sset\subsubsection{}\label{subsec_der-qcoh}
There are obvious forgetful functors:
$$
\iota_X:\cO_{\!X}\Mod_\qcoh\to\cO_{\!X}\Mod \qquad
R\iota_X:\sD(\cO_{\!X}\Mod_\qcoh)\to\sD(\cO_{\!X}\Mod)_\qcoh
$$
and we wish to exhibit right adjoints to these functors.
To this aim, suppose first that $X$ is affine; then we
may consider the functor:
$$
\sqcoh_X:\cO_{\!X}\Mod\to\cO_{\!X}\Mod_\qcoh
\qquad \cF\mapsto\cF(X)^\sim
$$
If $\cG$ is a quasi-coherent $\cO_{\!X}$-module, then
clearly $\sqcoh_X\cG\simeq\cG$; moreover, for any other
$\cO_{\!X}$-module $\cF$, there is a natural bijection:
$$
\Hom_{\cO_{\!X}}(\cG,\cF)\isom\Hom_{\cO_{\!X}(X)}(\cG(X),\cF(X)).
$$
It follows easily that $\sqcoh_X$ is the sought right adjoint.

\sset\subsubsection{}\label{subsec_cohereur-q-aff}
Slightly more generally, let $U$ be {\em quasi-affine\/},
{\em i.e.} a quasi-compact open subset of an affine scheme,
and choose a quasi-compact open immersion $j:U\to X$ into
an affine scheme $X$. In this case, we may define
$$
\sqcoh_U:\cO_{\!U}\Mod\to\cO_{\!X}\Mod\qquad
\cF\mapsto(\sqcoh_X j_*\cF)_{|U}.
$$
Since $j_*:\cO_{\!U}\Mod\to\cO_{\!X}\Mod$ is right adjoint to
$j^*:\cO_{\!X}\Mod\to\cO_{\!U}\Mod$, we have a natural isomorphism:
$\Hom_{\cO_{\!U}}(\cG,\cF)\isom\Hom_{\cO_{\!X}}(j_*\cG,j_*\cF)$
for every $\cO_{\!U}$-modules $\cG$ and $\cF$. Moreover,
if $\cG$ is quasi-coherent, the same holds for $j_*\cG$
(\cite[Ch.I, Cor.9.2.2]{EGAI}), whence a natural isomorphism
$\cG\simeq\sqcoh_U\cG$, by the foregoing discussion for
the affine case. Summing up, this shows that $\sqcoh_U$
is a right adjoint to $\iota_U$, and especially it is
independent, up to unique isomorphism, of the choice of $j$.

\sset\subsubsection{}\label{subsec_cohereur-gen}
Next, suppose that $X$ is quasi-compact and quasi-separated.
We choose a finite covering $\fU:=(U_i~|~i\in I)$ of $X$, consisting
of affine open subsets, and set $U:=\coprod_{i\in I}U_i$, the
$X$-scheme which is the disjoint union ({\em i.e.} categorical
coproduct) of the schemes $U_i$. We denote by $\fU_\bullet$ the
simplicial covering such that $\fU_n:=U\times_X\cdots\times_XU$,
the $(n+1)$-th power of $U$, with the face maps given by the
natural projections, and degeneracies induced by the
diagonal map $U\to U\times_XU$; let also $\pi_n:\fU_n\to X$
be the natural morphism, for every $n\in\N$.
Clearly we have $\pi_{n-1}\circ\partial_i=\pi_n$ for every
face morphism $\partial_i:\fU_n\to\fU_{n-1}$. The simplicial
scheme $\fU_\bullet$ (with the Zariski topology on each
scheme $\fU_n$) can also be regarded as a fibred topos
over the category $\Delta^{\!o}$ (notation of \cite[\S2.2]{Ga-Ra});
then the datum $\cO_{\!\fU_\bullet}:=(\cO_{\!\fU_n}~|~n\in\N)$
consisting of the structure sheaves on each $\fU_n$ and the
natural maps $\partial_i^*\cO_{\!\fU_{n-1}}\to\cO_{\!\fU_n}$ for
every $n>0$ and every $i=0,\dots,n$ (and similarly for the
degeneracy maps), defines a ring in the associated topos
$\rmTop(\fU_\bullet)$ (see \cite[\S3.3.15]{Ga-Ra}).
We denote by $\cO_{\!\fU_\bullet}\Mod$ the category of
$\cO_{\!\fU_\bullet}$-modules in the topos $\rmTop(\fU_\bullet)$.
The family $(\pi_n~|~n\in\N)$ induces a morphism of topoi
$$
\pi_\bullet:\rmTop(\fU_\bullet)\to s.X
$$
where $s.X$ is the topos $\rmTop(X_\bullet)$ associated with the
constant simplicial scheme $X_\bullet$ (with its Zariski topology)
such that $X_n:=X$ for every $n\in\N$ and such that all the face
and degeneracy maps are $\one_X$. Clearly the objects of $s.X$
are nothing else than the cosimplicial Zariski sheaves on $X$.
Especially, if we view a $\cO_{\!X}$-module $\cF$ as a constant
cosimplicial $\cO_{\!X}$-module, we may define the
{\em augmented cosimplicial \v{C}ech $\cO_{\!X}$-module\/}
\set\begin{equation}\label{eq_aug-cos-Cech}
\cF\to\cC^\bullet(\fU,\cF):=\pi_{\bullet*}\circ\pi_\bullet^*\cF
\end{equation}
associated with $\cF$ and the covering $\fU$; in every degree
$n\in\N$ this is defined by the rule :
$$
\cC^n(\fU,\cF):=\pi_{n*}\pi^*_n\cF
$$
and the coface operators $\partial^i$ are induced by the face
morphisms $\partial_i$ in the obvious way.

\begin{lemma}\label{lem_aspher}
{\em (i)}\ \ 
The augmented complex \eqref{eq_aug-cos-Cech} is aspherical
for every $\cO_{\!X}$-module $\cF$.
\begin{enumerate}
\addenu
\item
If $\cF$ is an injective $\cO_{\!X}$-module, then
\eqref{eq_aug-cos-Cech} is homotopically trivial.
\end{enumerate}
\end{lemma}
\begin{proof} (See also {\em e.g.} \cite[Th.5.2.1]{God}.)
The proof relies on the following alternative description
of the cosimplicial \v{C}ech $\cO_{\!X}$-module.
Consider the adjoint pair :
$$
\xymatrix{(\pi^*,\pi_*):\cO_{\!\fU}\Mod \ar@<.5ex>[r] &
\cO_{\!X}\Mod \ar@<.5ex>[l]}
$$
arising from our covering $\pi:\fU\to X$.
Let $(\top:=\pi_*\circ\pi^*,\eta,\mu)$ be the associated
triple (see \eqref{subsec_perp-adj}), $\cF$ any
$\cO_{\!X}$-module; we leave to the reader the verification
that the resulting augmented cosimplicial complex
$\cF\to\top^\bullet\cF$ -- as defined in
\eqref{subsec_simpl-cotriple} -- is none else than the
augmented \v{C}ech complex \eqref{eq_aug-cos-Cech}.
Thus, for every $\cO_{\!X}$-module $\cF$, the augmented
complex $\pi^*\cF\to\pi^*\cC^\bullet(\fU,\cF)$ is homotopically
trivial (proposition \ref{prop_triple-res});
since $\pi$ is a covering morphism, (i) follows.
Furthermore, by the same token, the augmented complex
$\pi_*\cG\to\cC^\bullet(\fU,\pi_*\cG)$ is homotopically trivial
for every $\cO_{\!\fU}$-module $\cG$; especially we may take
$\cG:=\pi^*\cI$, where $\cI$ is an injective $\cO_{\!X}$-module.
On the other hand, when $\cI$ is injective, the unit of adjunction
$\cI\to\top\cI$ is split injective; hence the augmented complex
$\cI\to\cC^\bullet(\fU,\cI)$ is a direct summand of the
homotopically trivial complex $\top\cI\to\cC^\bullet(\fU,\top\cI)$,
and (ii) follows.
\end{proof}

Notice now that the schemes $\fU_n$ are quasi-affine for every
$n\in\N$, hence the functors $\sqcoh_{\fU_n}$ are well defined
as in \eqref{subsec_cohereur-q-aff}, and indeed, the rule :
$(\cF_n~|~n\in\N)\mapsto(\sqcoh_{\fU_n}\cF_n~|~n\in\N)$ yields
a functor :
$$
\sqcoh_{\fU_\bullet}:\cO_{\!\fU_\bullet}\Mod\to\cO_{\!\fU_\bullet}\Mod.
$$
This suggests to introduce a {\em quasi-coherent cosimplicial
\v{C}ech $\cO_{\!X}$-module} :
$$
\qcC^\bullet(\fU,\cF):=
\pi_{\bullet*}\circ\sqcoh_{\fU_\bullet}\circ\pi_\bullet^*\cF
$$
for every $\cO_{\!X}$-module $\cF$ (regarded as a constant cosimplicial
module in the usual way). More plainly, this is the cosimplicial
$\cO_{\!X}$-module such that :
$$
\qcC^n(\fU,\cF):=\pi_{n*}\circ\sqcoh_{\fU_n}\circ\pi^*_n\cF
\qquad\text{for every $n\in\N$}.
$$
According to \cite[Ch.I, Cor.9.2.2]{EGAI}, the $\cO_{\!X}$-modules
$\qcC^n(\fU,\cF)$ are quasi-coherent for all $n\in\N$.
Finally, we define the functor :
$$
\sqcoh_X:\cO_{\!X}\Mod\to\cO_{\!X}\Mod_\qcoh\qquad\cF\mapsto
\xymatrix{\Equal(\qcC^0(\fU,\cF)\ar@<.5ex>[r]^-{\partial^0}
\ar@<-.5ex>[r]_-{\partial^1} &\qcC^1(\fU,\cF))}
$$

\begin{proposition}\label{prop_cohereur}
{\em (i)}\ \ 
In the situation of \eqref{subsec_cohereur-gen},
the functor $\sqcoh_X$ is right adjoint to $\iota_X$.
\begin{enumerate}
\addenu
\item
Let $Y$ be any other quasi-compact and quasi-separated scheme,
$f:X\to Y$ any morphism. Then the induced diagram of\/ functors:
$$
\xymatrix{
\cO_{\!X}\Mod \ar[r]^-{f_*} \ar[d]_{\sqcoh_X} &
\cO_Y\Mod \ar[d]^{\sqcoh_Y} \\
\cO_{\!X}\Mod_\qcoh \ar[r]^-{f_*} & \cO_Y\Mod_\qcoh
}$$
commutes up to a natural isomorphism of functors.
\end{enumerate}
\end{proposition}
\begin{proof} For every $n\in\N$ and every $\cO_{\!\fU_n}$-module
$\cH$, the counit of the adjunction yields a natural map of
$\cO_{\!\fU_n}$-modules: $\sqcoh_{\fU_n}\cH\to\cH$. Taking $\cH$
to be $\pi_n^*\cF$ on $\fU_n$ (for a given $\cO_{\!X}$-module $\cF$),
these maps assemble to a morphism of cosimplicial $\cO_{\!X}$-modules :
\set\begin{equation}\label{eq_qcC-versus-cC}
\qcC^\bullet(\fU,\cF)\to\cC^\bullet(\fU,\cF)
\end{equation}
and it is clear that \eqref{eq_qcC-versus-cC} is an isomorphism
whenever $\cF$ is quasi-coherent. Let now $\phi:\cG\to\cF$ be
a map of $\cO_{\!X}$-modules, with $\cG$ quasi-coherent; after
applying the natural transformation \eqref{eq_qcC-versus-cC}
and forming equalizers, we obtain a commutative diagram :
$$
\xymatrix{
\sqcoh_X\cG \ar[r]^-\sim \ar[d]_-{\sqcoh_X\phi} & \cG \ar[d]^-\phi \\
 \sqcoh_X\cF \ar[r] & \cF
}$$
from which we see that the rule: $\phi\mapsto\sqcoh_X\phi$
establishes a natural injection:
\set\begin{equation}\label{eq_nat-inject}
\Hom_{\cO_{\!X}}(\cG,\cF)\to\Hom_{\cO_{\!X}}(\cG,\sqcoh_X\cF)
\end{equation}
and since $\Hom_{\cO_{\!X}}(\cG,\eqref{eq_qcC-versus-cC})$ is an
isomorphism of cosimplicial $\cO_{\!X}$-modules,
\eqref{eq_nat-inject} is actually a bijection, whence (i).

(ii) is obvious, since both $\sqcoh_Y\circ f_*$ and
$f_*\circ\sqcoh_X$ are right adjoint to
$f^*\circ\iota_Y=\iota_X\circ f^*$.
\end{proof}

\sset\subsubsection{}
In the situation of \eqref{subsec_cohereur-gen},
the functor $\sqcoh_X$ is left exact, since it is a right
adjoint, hence it gives rise to a left derived functor :
$$
R\sqcoh_X:\sD^+(\cO_{\!X}\Mod)\to\sD^+(\cO_{\!X}\Mod_\qcoh).
$$

\begin{proposition}\label{prop_der-cohereur}
Let $X$ be a quasi-compact and quasi-separated scheme. Then :
\begin{enumerate}
\item
$R\sqcoh_X$ is right adjoint to $R\iota_X$.
\item
Suppose moreover, that $X$ is {\em semi-separated}, {\em i.e.}
such that the intersection of any two affine open subsets of $X$,
is still affine. Then the unit of the adjunction $(R\iota_X,R\sqcoh_X)$
is an isomorphism of functors.
\end{enumerate}
\end{proposition}
\begin{proof} (i): The exactness of the functor $\iota_X$
implies that $\sqcoh_X$ preserves injectives; the assertion
is a formal consequence : the details shall be left to the reader.

(ii): Let $\cF^\bullet$ be any complex of quasi-coherent
$\cO_{\!X}$-modules; we have to show that the natural map
$\cF^\bullet\to R\sqcoh_X\cF^\bullet$ is an isomorphism.
Using a Cartan-Eilenberg resolution
$\cF^\bullet\isom\cI^{\bullet\bullet}$ we deduce a spectral
sequence (\cite[\S5.7]{We})
$$
E_1^{pq}:=R^p\sqcoh_X H^q\cF^\bullet\Rightarrow
R^{p+q}\sqcoh_X\cF^\bullet
$$
which easily reduces to checking the assertion for
the cohomology of $\cF^\bullet$, so we may assume from
start that $\cF^\bullet$ is a single $\cO_{\!X}$-module
placed in degree zero. Let us then choose an injective
resolution $\cF\isom\cI^\bullet$ (that is, in the category
of all $\cO_{\!X}$-modules); we have to show that
$H^p\sqcoh_X\cI^\bullet=0$ for $p>0$. We deal first
with the following special case :

\begin{claim}\label{cl_aff-case}
Assertion (i) holds if $X$ is affine.
\end{claim}
\begin{pfclaim} Indeed, in this case, the chosen injective
resolution of $\cF$ yields a long exact sequence
$0\to\cF(X)\to\cI^\bullet(X)$, and therefore a resolution
$\sqcoh_X\cF:=\cF(X)^\sim\to\sqcoh_X\cI^\bullet:=\cI^\bullet(X)^\sim$.
\end{pfclaim}

For the general case, we choose any affine covering $\fU$ of $X$
and we consider the cochain complex of cosimplicial complexes
$\qcC^\bullet(\fU,\cI^\bullet)$.

\begin{claim}\label{cl_contract} For any injective
$\cO_{\!X}$-module $\cI$, the augmented complex:
$$
\sqcoh_X\cI\to\qcC^\bullet(\fU,\cI)
$$
is homotopically trivial.
\end{claim}
\begin{pfclaim} It follows easily from proposition
\ref{prop_cohereur}(ii) that
$$
\qcC^\bullet(\fU,\cF)\simeq\sqcoh_X(\cC^\bullet(\fU,\cF))
$$
for any $\cO_{\!X}$-module $\cF$. Then the claim follows
from lemma \ref{lem_aspher}(ii).
\end{pfclaim}

We have a spectral sequence :
$$
E_1^{pq}:=H^p\qcC^\bullet(\fU,\cI^q)\Rightarrow
\Tot^{p+q}(\qcC^\bullet(\fU,\cI^\bullet))
$$
and it follows from claim \ref{cl_contract} that
$E^{pq}_1=0$ whenever $p>0$, and $E^{0q}=\sqcoh_X\cI^q$ for
every $q\in\N$, so the spectral sequence $E^{\bullet\bullet}$
degenerates, and we deduce a quasi-isomorphism
\set\begin{equation}\label{eq_Tot-is-inj}
\sqcoh_X\cI^\bullet\isom
\Tot^\bullet(\qcC^\bullet(\fU,\cI^\bullet)).
\end{equation}
On the other hand, for fixed $q\in\N$, we have
$\qcC^q(\fU,\cI^\bullet)=\pi_{n*}\sqcoh_{\fU_n}\pi^*_n\cI^\bullet$;
since $X$ is semi-separated, $\fU_n$ is affine, so
the complex $\sqcoh_{\fU_n}\pi^*_n\cI^\bullet$ is a resolution
of $\sqcoh_{\fU_n}\pi^*_n\cF=\pi^*_n\cF$, by claim \ref{cl_aff-case}.
Furthermore, $\pi_n:\fU_n\to X$ is an affine morphism,
so $\qcC^q(\fU,\cI^\bullet)$ is a resolution of $\pi_{n*}\pi^*_n\cF$.
Summing up, we see that
$\Tot^\bullet(\qcC^\bullet(\fU,\cI^\bullet))$ is
quasi-isomorphic to $\cC^\bullet(\fU,\cF)$, which is
a resolution of $\cF$, by lemma \ref{lem_aspher}(i).
Combining with \eqref{eq_Tot-is-inj}, we deduce (ii).
\end{proof}

\begin{theorem}\label{th_cohereur}
Let $X$ be a quasi-compact semi-separated scheme.
The forgetful functor
$$
R\iota_X:\sD^+(\cO_{\!X}\Mod_\qcoh)\to\sD^+(\cO_{\!X}\Mod)_\qcoh
$$
is an equivalence of categories, whose quasi-inverse is the
restriction of $R\sqcoh_X$.
\end{theorem}
\begin{proof} By proposition \ref{prop_der-cohereur}(ii) we know
already that the composition $R\sqcoh_X\circ R\iota_X$ is a
self-equivalence of $\sD^+(\cO_{\!X}\Mod_\qcoh)$. For every
complex $\cF^\bullet$ in $\sD^+(\cO_{\!X}\Mod)$, the counit
of adjunction $\eps_{\cF^\bullet}:
R\iota_X\circ R\sqcoh_X\cF^\bullet\to\cF^\bullet$
can be described as follows. Pick an injective resolution
$\alpha:\cF^\bullet\to\cI^\bullet$; then $\eps_{\cF^\bullet}$
is defined by the diagram :
$$
\xymatrix{
\sqcoh_X\cI^\bullet \ar[r]^-{\eps^\bullet} & \cI^\bullet &
\cF^\bullet \ar[l]_-\alpha
}$$
where, for each $n\in\N$, the map $\eps^n:\sqcoh_X\cI^n\to\cI^n$
is the counit of the adjoint pair $(\iota_X,\sqcoh_X)$. It
suffices then to show that $\eps^\bullet$ is a quasi-isomorphism,
when $\cI^\bullet$ is an object of $\sD^+(\cO_{\!X}\Mod)_\qcoh$.
To this aim, we may further choose $\cI^\bullet$ of the form
$\Tot^\bullet(\cI^{\bullet\bullet})$, where
$\cI^{\bullet\bullet}$ is a Cartan-Eilenberg resolution
of $\cF^\bullet$ (see \cite[\S5.7]{We}). We then deduce
a spectral sequence :
$$
E^{pq}_2:=R^p\sqcoh_XH^q\cF^\bullet\Rightarrow
R^{p+q}\sqcoh_X\cF^\bullet.
$$
The double complex $\cI^{\bullet\bullet}$ also gives
rise to a similar spectral sequence $F_2^{pq}$, and
clearly $F^{pq}_2=0$ whenever $p>0$, and $F^{0q}=H^q\cF^\bullet$.
Furthermore, the counit of adjunction
$\eps^{\bullet\bullet}:\sqcoh_X\cI^{\bullet\bullet}\to
\cI^{\bullet\bullet}$ induces a morphism of spectral
sequences $\omega^{pq}:E^{pq}_2\to F^{pq}_2$.
Consequently, in order to prove that the $\eps_{\cF^\bullet}$
is a quasi-isomorphism, it suffices to show that $\omega^{pq}$
is an isomorphism for every $p,q\in\N$. This comes down
to the assertion that $R^p\sqcoh_X\cG=0$ for every
quasi-coherent $\cO_{\!X}$-module $\cG$, and every  $p>0$.
However, we have $\cG=R\iota_X\cG$, so that
$R\sqcoh_X\cG=R\sqcoh_X\circ R\iota_X\cG$, and then
the contention follows from proposition \ref{prop_der-cohereur}(ii).
\end{proof}

\begin{lemma}\label{lem_extend-cohs}
Let $X$ be a quasi-compact and quasi-separated scheme,
$U$ a quasi-compact open subset of $X$, $\cH$ a quasi-coherent
$\cO_{\!X}$-module, $\cG$ a finitely presented quasi-coherent
$\cO_{\!U}$-module, and $\phi:\cG\to\cH_{|U}$ an $\cO_{\!U}$-linear
map. Then:
\begin{enumerate}
\item
There exist a finitely presented quasi-coherent
$\cO_{\!X}$-module $\cF$ on $X$, and a $\cO_{\!X}$-linear map
$\psi:\cF\to\cH$, such that $\cF_{|U}=\cG$ and $\psi_{|U}=\phi$.
\item
Especially, every finitely presented quasi-coherent
$\cO_{\!U}$-module extends to a finitely presented quasi-coherent
$\cO_{\!X}$-module.
\end{enumerate}
\end{lemma}
\begin{proof} (i): Let $(V_i~|~i=1,\dots,n)$ be a finite
affine open covering of $X$. For every $i=0,\dots,n$, let us
set $U_i:=U\cup V_1\cup\cdots\cup V_i$; we construct, by
induction on $i$, a family of finitely presented
quasi-coherent $\cO_{\!U_i}$-modules $\cF_i$, and morphisms
$\psi_i:\cF_i\to\cH_{|U_i}$ such that $\cF_{i+1|U_i}=\cF_i$
and $\psi_{i+1|U_{i+1}}=\psi_i$ for every $i<n$. For $i=0$
we have $U_0=U$, and we set $\cF_0:=\cG$, $\psi_0:=\phi$.
Suppose that $\cF_i$ and $\psi_i$ have already been given.
Since $X$ is quasi-separated, the same holds for $U_{i+1}$,
and the immersion $j:U_i\to U_{i+1}$ is quasi-compact;
it follows that $j_*\cF_i$ and $j_*\cH_{|U_i}$ are quasi-coherent
$\cO_{\!U_i}$-modules (\cite[Ch.I, Prop.9.4.2(i)]{EGAI}).
We let
$$
\cM:=j_*\cF_i\times_{j_*\cH_{|U_i}}\cH_{|U_{i+1}}.
$$
Then $\cM$ is a quasi-coherent $\cO_{\!U_{i+1}}$-module
admitting a map $\cM\to\cH_{|U_{i+1}}$, and such that
$\cM_{|U_i}=\cF_i$. We can then find a filtered family
of quasi-coherent $\cO_{V_{i+1}}$-modules of finite
presentation $(\cM_\lambda~|~\lambda\in\Lambda)$, whose
colimit is $\cM_{|V_{i+1}}$. Since $\cF_i$ is finitely
presented and $U_i\cap V_{i+1}$ is quasi-compact, there
exists $\lambda\in\Lambda$ such that the induced morphism
$\beta:\cM_{\lambda|U_i\cap V_{i+1}}\to\cF_{i|U_i\cap V_{i+1}}$
is an isomorphism. We can thus define $\cF_{i+1}$ by
gluing $\cF_i$ and $\cM_\lambda$ along $\beta$; likewise,
$\psi_i$ and the induced map
$\cM_\lambda\to\cM_{|V_{i+1}}\to\cH_{|V_{i+1}}$ glue
to a map $\psi_{i+1}$ as required.
Clearly the pair $(\cF:=\cF_n,\psi:=\psi_n)$ is the sought
extension of $(\cG,\phi)$.

(ii): Let $0_X$ be the final object in the category of
$\cO_{\!X}$-modules; to extend a finitely presented quasi-coherent
$\cO_{\!U}$-module $\cG$, it suffices to apply (i) to the
unique map $\cG\to 0_{X|U}$.
\end{proof}

For ease of reference, we point out the following simple
consequence of lemma \ref{lem_extend-cohs}.

\begin{corollary}\label{cor_ease-of-ref}
Let $U$ be a quasi-affine scheme, $\cE$ a locally free
$\cO_{\!U}$-module of finite type. Then we may find integers
$n,m\in\N$ and a left exact sequence of $\cO_{\!U}$-modules :
$$
0\to\cE\to\cO_{\!U}^{\oplus m}\to\cO_{\!U}^{\oplus n}.
$$
\end{corollary}
\begin{proof} Set $\cE^\vee:=\cHom_{\cO_{\!U}}(\cE,\cO_{\!U})$,
and notice that $\cE^\vee$ is a locally free $\cO_{\!U}$-module
of finite type, especially it is quasi-coherent and of finite
presentation.
By assumption, $U$ is a quasi-compact open subset of an affine
scheme $X$, hence $\cE^\vee$ extends to a quasi-coherent
$\cO_{\!X}$-module $\cF$ of finite presentation (lemma
\ref{lem_extend-cohs}(ii)); we have $\cF=F^\sim$, for some
finitely presented $A$-module $F$; we choose a
presentation of $F$ as the cokernel of an $A$-linear map
$A^{\oplus n}\to A^{\oplus m}$, for some $m,n\in\N$, whence
a presentation of $\cE^\vee$ as the cokernel of a morphism
$\cO_{\!U}^{\oplus n}\to\cO_{\!U}^{\oplus m}$, and after
dualizing again, we get the sought left exact sequence.
\end{proof}

\sset\subsubsection{}\label{subsec_fp-approx}
Let $X$ be a scheme; for any quasi-coherent $\cO_{\!X}$-module
$\cF$, we consider the full subcategory $C/\cF$ of the category
$\cO_{\!X}\Mod_\qcoh/\cF$ whose objects are all the maps
$\phi:\cG\to\cF$ such that $\cG$ is a finitely presented
$\cO_{\!X}$-module (notation of \eqref{subsec_slice-cat}).
We denote by
$$
\iota_\cF:C/\cF\to\cO_{\!X}\Mod_\qcoh
\qquad
(\cG\to\cF)\mapsto\cG
$$
the restriction of the source functor. The first
observation is:

\begin{lemma}\label{lem_filtered-cat}
$C/\cF$ is a filtered category.
\end{lemma}
\begin{proof} Since $C/\cF$ always admits the initial
object $0\to\cF$, its set of objects is non-empty.
Thus, according to remark \ref{rem_cofinal}(i),
it suffices to check that $C/\cF$ is directed and
satisfies the coequalizing condition. However, say
that $\phi_i:\cG_i\to\cF$ ($i=1,2$) are any two objects of
$C/\cF$; then obviously there is a unique morphism
$\cG_1\oplus\cG_2\to\cF$ whose restriction to $\cG_i$
agrees with $\phi_i$ for $i=1,2$.

Next, suppose that $\beta_1,\beta_2:\cG\to\cG'$ are two
morphisms of finitely presented $\cO_{\!X}$-modules
and $\psi:\cG'\to\cF$ an object of $C/\cF$ such that
$\psi\circ\beta_1=\psi\circ\beta_2$. Set
$\cG'':=\Coker\,(\beta_1-\beta_2)$; it is easily seen
that $\cG''$ is a finitely presented $\cO_{\!X}$-module
and $\psi$ factors through a unique morphism $\cG''\to\cF$,
whence the lemma.
\end{proof}

\begin{proposition}\label{prop_fp-approx}
With the notation of \eqref{subsec_fp-approx}, suppose that
$X$ is quasi-compact and quasi-separated. Then, for any
quasi-coherent $\cO_{\!X}$-module $\cF$, the induced map
$$
\colim_{C/\cF}\iota_\cF\to\cF
$$
is an isomorphism.
\end{proposition}
\begin{proof} Let $\cF'$ be such colimit; clearly there is
a natural map $\cF'\to\cF$, and we have to show that it is
an isomorphism. To this aim, we can check on the stalk
over the points $x\in X$, hence we come down to showing :

\begin{claim} Let $\cF$ be any quasi-coherent $\cO_{\!X}$-module.
Then :
\begin{enumerate}
\item
For every $s\in\cF_x$ there exist a finitely presented
quasi-coherent $\cO_{\!X}$-module $\cG$, a morphism
$\phi:\cG\to\cF$ and $t\in\cG_x$ such that $\phi_x(t)=s$.
\item
For every map $\phi:\cG\to\cF$ with $\cG$ finitely
presented and quasi-coherent, and every $s\in\Ker\,\phi_x$,
there exists a commutative diagram of quasi-coherent
$\cO_{\!X}$-modules :
$$
\xymatrix{
\cG \ar[dr]^\phi \ar[d]_\psi \\
\cH \ar[r] & \cF
}$$
with $\cH$ finitely presented and $\psi_x(s)=0$.
\end{enumerate}
\end{claim}
\begin{pfclaim}[](i): Let $U\subset X$ be an open
subset such that $s$ extends to a section $s_U\in\cF(U)$;
we deduce a map $\phi_U:\cO_{\!U}\to\cF_{|U}$ by the rule
$a\mapsto a\cdot s_U$ for every $a\in\cO_{\!U}(U)$.
In view of lemma \ref{lem_extend-cohs}(i), the pair
$(\phi_U,\cO_{\!U})$ extends to a pair $(\phi,\cG)$
with the sought properties.

(ii): We apply (i) to the quasi-coherent $\cO_{\!X}$-module
$\Ker\,\phi$, to find a finitely presented quasi-coherent
$\cO_{\!X}$-module $\cG'$, a morphism $\beta:\cG'\to\Ker\,\phi$
and $t\in\cH_x$ such that $\beta_x(t)=s$. Then we let
$\psi:\cG\to\cH:=\Coker(\cG'\xrightarrow{\beta}\Ker\,\phi\to\cG)$
be the natural map. By construction, $\cH$ is finitely presented,
$\psi_x(s)=0$ and clearly $\phi$ factors through $\psi$.
\end{pfclaim}
\end{proof}

\sset\subsubsection{}\label{subsec_coher-then-noether}
Let $X$ be a coherent scheme, $U\subset X$ a quasi-compact
open subset, $a,b\in\Z$ any two integers with $a\leq b$, and
$\cH^\bullet$ any object of $\sD^{[a,b]}(\cO_{\!X}\Mod)_\qcoh$ such
that $\cH^\bullet_{|U}$ lies in $\sD^{[a,b]}(\cO_{\!U}\Mod)_\coh$.
Let also $\cF_i\subset H^i\cH^\bullet$ be a quasi-coherent
$\cO_{\!X}$-submodule of finite type, for every $i=a,\dots,b$.

\begin{proposition}\label{prop_coher-then-noether}
In the situation of \eqref{subsec_coher-then-noether},
the following holds :
\begin{enumerate}
\item
There exists an object $\cG^\bullet$ of
$\sD^{[a,b]}(\cO_{\!X}\Mod)_\coh$, with a morphism
$\phi^\bullet:\cG^\bullet\to\cH^\bullet$ in
$\sD(\cO_{\!X}\Mod)_\qcoh$ such that :
\begin{enumerate}
\item
$\phi^\bullet_{|U}$ is an isomorphism.
\item
$\cF_i\subset\Img\,H^i\phi^\bullet$\ \ for every $i=a,\dots,b$.
\end{enumerate}
\item
If $X$ is noetherian, we can find $\phi^\bullet$ as in
{\em (i)} such that $H^i\phi^\bullet$ is a monomorphism,
for every $i=a,\dots,b$.
\end{enumerate}
\end{proposition}
\begin{proof} To start out, for any quasi-coherent
$\cO_{\!X}$-module $\cF$ such that $\cF_{|U}$ is a
coherent $\cO_{\!U}$-module, define the category
$C/\cF$ as in \eqref{subsec_fp-approx}, and consider
the full subcategory $C'/\cF$ (resp. $C''/\cF$) of
$C/\cF$ whose objects are all the maps $\phi:\cG\to\cF$
such that $\phi_{|U}$ is an isomorphism (resp. an epimorphism).

\begin{claim} $C'/\cF$ and $C''/\cF$ are filtered
categories.
\end{claim}
\begin{pfclaim} The proof of lemma \ref{lem_filtered-cat}
applies {\em verbatim} to show that both $C'/\cF$ and
$C''/\cF$ satisfy the coequalizing condition of definition
\ref{def_filtered-cols}(v), and it also proves that $C''/\cF$
is directed. Moreover, lemma \ref{lem_extend-cohs}(i) says
that the sets of objects of these two categories are non-empty,
which already implies that $C''/\cF$ is filtered, by
remark \ref{rem_cofinal}(i). By the same token, it remains
only to check that $C'/\cF$ satisfies the coequalizing
condition. Hence, let $\phi_i:\cG_i\to\cF$ ($i=1,2$) be any
two objects of $C'/\cF$, and set $\cG':=\cG_1\times_\cF\cG_2$;
then $\cG'$ is a quasi-coherent $\cO_{\!X}$-module, and
$\cG'_{|U}$ is isomorphic to $\cF_{|U}$, especially it is
coherent. Then pick any object $\psi:\cG''\to\cG'$ of the
category $C''/\cG'$, and denote by $\psi_i:\cG''\to\cG_i$
($i=1,2$) the two morphisms induced by $\psi$; by
construction we have $\phi_1\circ\psi_1=\phi_2\circ\psi_2$.
Hence, set $\cG:=\cG_1\amalg_{\cG''}\cG_2$; we get a
unique morphism $\phi:\cG\to\cF$ such that both $\phi_1$
and $\phi_2$ factor through $\phi$ and $\phi_{|U}$ is
an isomorphism. Since $\cG$ is a finitely presented
$\cO_{\!X}$-module, we are done.
\end{pfclaim}

\begin{claim}\label{cl_fp-approx}
With the foregoing notation, the following holds :
\begin{enumerate}
\item
$C'/\cF$ is a cofinal subcategory of $C/\cF$.
\item
For every object $K^\bullet$ of $\sD^b(\cO_{\!X}\Mod)_\coh$
and every $c\in\Z$, the natural map
\set\begin{equation}\label{eq_filtered-ab-group}
\colim_{C/\cF}\Hom_{\sD(\cO_{\!X}\Mod)}(K^\bullet,\iota_\cF[c])
\to\Hom_{\sD(\cO_{\!X}\Mod)}(K^\bullet,\cF[c])
\end{equation}
is an isomorphism (notation of \eqref{subsec_fp-approx}).
\end{enumerate}
\end{claim}
\begin{pfclaim}(i): Since $U$ is quasi-compact, it is easily
seen that $C''/\cF$ is cofinal in $C/\cF$, so we are reduced to
checking that $C'/\cF$ is cofinal in $C''/\cF$. Hence, let
$\psi:\cG'\to\cF$ be any object of $C''/\cF$; since $\cO_{\!X}$
is coherent, $\cK:=\Ker\,\psi_{|U}$ is a coherent $\cO_{\!U}$-module,
therefore we may extend $\cK$ to a coherent $\cO_{\!X}$-module
$\cK'$ and the identity map of $\cK$ to a morphism
$\psi':\cK'\to\Ker\,\psi$ of $\cO_{\!X}$-modules (lemma
\ref{lem_extend-cohs}(i)). Set $\cG:=\cG'/\psi'(\cK')$;
the induced map $\cG\to\cF$ is an object of $C'/\cF$, whence
the contention.

(ii): We have two spectral sequences
$$
\begin{aligned}
E^{pq}_2 :=\, & \colim_{C/\cF}
R^p\Gamma\circ R^q\cHom^\bullet_{\cO_{\!X}}(K^\bullet,\iota_\cF[c])
\Rightarrow\colim_{C/\cF}
\Hom_{\sD(\cO_{\!X}\Mod)}(K^\bullet,\iota_\cF[c+p+q]) \\
F^{pq}_2 :=\, &
R^p\Gamma\circ R^q\cHom^\bullet_{\cO_{\!X}}(K^\bullet,\cF[c])
\Rightarrow\Hom_{\sD(\cO_{\!X}\Mod)}(K^\bullet,\cF[c+p+q])
\end{aligned}
$$
and a morphism of spectral sequences
$E^{\bullet\bullet}_\bullet\to F^{\bullet\bullet}_\bullet$,
such that \eqref{eq_filtered-ab-group} is a morphism
of filtered abelian groups, for the two finite filtrations
induced by these spectral sequences on their abutments.
Since the functors $R^p\Gamma$ commute with filtered colimits,
we deduce that it suffices to show that the natural morphism
$$
\colim_{C/\cF}R^q\cHom^\bullet_{\cO_{\!X}}(K^\bullet,\iota_\cF[c]))
\to R^q\cHom^\bullet_{\cO_{\!X}}(K^\bullet,\cF[c])
$$
is an isomorphism for every $q\in\Z$. Then, a standard
{\em d\'evissage\/} argument further reduces to the case
where $K^\bullet$ is concentrated in a single degree,
so we come down to checking that the functor
$\cF\mapsto R^q\cHom^\bullet_{\cO_{\!X}}(\cG,\cF)$ commutes
with filtered colimits of quasi-coherent $\cO_{\!X}$-modules,
for every coherent $\cO_{\!X}$-module $\cG$ and every
$q\in\Z$. To this aim we may assume that $X$ is affine,
in which case -- since $\cO_{\!X}$ is coherent -- $\cG$
admits a resolution $\cL^\bullet$ consisting of
free $\cO_{\!X}$-modules of finite rank; the latter are
acyclic for the functor $\cHom_{\cO_{\!X}}$, so we are
left with the assertion that the functor
$\cF\mapsto H^q\cHom_{\cO_{\!X}}(\cL^\bullet,\cF)$
commutes with filtered colimits, which is clear.
\end{pfclaim}

\begin{claim}\label{cl_a-is-b}
The proposition holds if $a=b$.
\end{claim}
\begin{pfclaim} From proposition \ref{prop_fp-approx} and
claim \ref{cl_fp-approx}(i) we deduce that, for every
$x\in X$ we may find an object $\phi:\cG\to H^a\cH$
such that $\cF_{a,x}\subset\Img\,\phi_x$, and therefore
$\cF_{a|U_x}\subset\Img\,\phi_{|U_x}$ for some open
neighborhood $U_x$ of $x$ in $X$. Since $X$ is quasi-compact
and $C/\cF$ is filtered, assertion (i) follows easily.
Next, suppose that a morphism $\phi:\cG\to H^a\cH$ has
already been exhibited that fulfills conditions (i.a)
and (i.b), and set $\cG':=\Img\,\phi$; to show assertion
(ii), it suffices to remark that if $X$ is noetherian,
$\cG'$ is still a coherent $\cO_{\!X}$-module.
\end{pfclaim}

We proceed now by induction on $b-a$, the case where $a=b$
having been covered by claim \ref{cl_a-is-b}.
Thus, suppose that $a<b$, and that the proposition is already
known for every object of $\sD^{[a+1,b]}(\cO_{\!X}\Mod)_\qcoh$
fulfilling the stated condition. We have a distinguished
triangle
$$
\cF[-a]\to\cH^\bullet\to\tau^{\geq a+1}\cH^\bullet
\xrightarrow{\ \alpha^\bullet\ }\cF[1-a]
\qquad
\text{with $\cF:=H^a\cH^\bullet$}
$$
in $\sD(\cO_{\!X}\Mod)_\qcoh$ (where $\tau^{\geq a+1}$ denotes
the truncation functor : see remark \ref{rem_derived-cat}(iv)).
By inductive assumption, we may find an object $\cG'{}^\bullet$
of $\sD^{[a+1,b]}(\cO_{\!X}\Mod)_\coh$ and a morphism
$\phi'{}^\bullet:\cG'{}^\bullet\to\tau^{\geq a+1}\cH^\bullet$ in
$\sD(\cO_{\!X}\Mod)_\qcoh$ which restricts to an isomorphism
on $U$ and such that $\cF_i\subset\Img\,H^i\phi'{}^\bullet$
for every $i=a+1,\dots,b$. There follows a morphism of
distinguished triangles :
$$
\xymatrix{
\cF[-a] \ar[r] \ddouble &
\cH'{}^\bullet \ar[r] \ar[d]^{\beta^\bullet} &
\cG'{}^\bullet \ar[r] \ar[d]^{\phi'{}^\bullet} &
\cF[1-a] \ddouble \\
\cF[-a] \ar[r] & \cH^\bullet \ar[r] &
\tau^{\geq a+1}\cH^\bullet \ar[r]^-{\alpha^\bullet} & \cF[1-a]
}$$
for some object $\cH'{}^\bullet$ of $\sD^{[a,b]}(\cO_{\!X}\Mod)_\qcoh$,
and clearly $\beta^\bullet$ restricts to an isomorphism on $U$
and $H^a\beta^\bullet$ is an isomorphism. We may then replace
$\cH^\bullet$ by $\cH'{}^\bullet$, and assume that
$\tau^{\geq a+1}\cH^\bullet$ lies in $\sD^{[a+1,b]}(\cO_{\!X}\Mod)_\coh$.

By claim \ref{cl_fp-approx}, we may find an object
$\psi:\cG\to\cF$ of $C'/\cF$ such that
$\alpha^\bullet$ factors through $\psi[1-a]$ and
a morphism $\omega^\bullet:\tau^{\geq a+1}\cH^\bullet\to\cG[1-a]$
in $\sD(\cO_{\!X}\Mod)_\coh$. Moreover, by claim
\ref{cl_a-is-b} we may assume that $\cF_a\subset\Img\,\psi$,
and also that $\psi$ is a monomorphism, if $X$ is
noetherian. Then $\Cone\,\omega^\bullet[-1]$ lies in
$\sD^{[a,b]}(\cO_{\!X}\Mod)_\coh$, and the induced
morphism $\Cone\,\omega^\bullet[-1]\to\cH^\bullet$ will do.
\end{proof}

\begin{corollary}\label{cor_fp-approx}
Let $X$ be a coherent scheme, $U\subset X$ a quasi-compact
open subset. Then the induced functor
$$
\sD^{[a,b]}(\cO_{\!X}\Mod)_\coh\to\sD^{[a,b]}(\cO_{\!U}\Mod)_\coh
$$
is essentially surjective, for every $a,b\in\Z$ with $a\leq b$.
\end{corollary}
\begin{proof} Let $j:U\to X$ be the open immersion,
$\cF^\bullet$ an object of $\sD^{[a,b]}(\cO_{\!U}\Mod)_\coh$.
By \cite[Ch.III, Prop.1.4.10, Cor.1.4.12]{EGAIII} and
a standard spectral sequence argument, it is easily
seen that
$$
\cH^\bullet:=\tau^{\geq a}\circ\tau^{\leq b}Rj_*\cF^\bullet
$$
lies in $\sD^{[a,b]}(\cO_{\!X}\Mod)_\qcoh$ and
$\cH^\bullet_{|U}=\cF^\bullet$. Then the assertion follows
immediately from proposition \ref{prop_coher-then-noether}(i).
\end{proof}

\subsection{Depth and cohomology with supports}\label{subsec_supp}
This section introduces and studies local cohomology and the
closely related notion of depth, in the context of arbitrary
schemes (whereas the usual references restrict to the case of
locally noetherian schemes).

\sset\subsubsection{}\label{subsec_Gamma_Z}
To begin with, let $(X,\cA)$ be any ringed topological space,
$i:Z\to X$ a closed immersion, set $U:=X\!\setminus\!Z$ and
denote by $j:U\to X$ the resulting open immersion. One defines
the functor
$$
\underline\Gamma_Z:\cA\Mod\to\cA\Mod \qquad
\cF\mapsto\Ker\,(\cF\to j_*j^*\cF)
$$
as well as its composition with the global section functor
$$
\Gamma_{\!\!Z}:\cA\Mod\to\cA(X)\Mod \qquad
\cF\mapsto\Gamma(X,\underline\Gamma_Z\cF).
$$
It is clear that $\underline\Gamma_Z$ and $\Gamma_{\!Z}$ are
left exact functors, hence they give rise to right derived
functors
$$
R\underline\Gamma_Z:\sD^+(\cA\Mod)\to\sD^+(\cA\Mod)
\qquad
R\Gamma_{\!\!Z}:\sD^+(\cA\Mod)\to\sD^+(\cA(X)\Mod)
$$
such that the natural morphism
$\underline\Gamma_Z\cF\to R^0\underline\Gamma_Z\cF$
is an isomorphism, for every $\cA$-module $\cF$
(and likewise for $R^0\Gamma_{\!Z}$).
Moreover, suppose that $\cF\to\cI^\bullet$ is a
resolution of $\cF$ by a complex of injective
$\cA$-modules; since injective sheaves are flabby
(lemma \ref{lem_franziska}(v)), we obtain a short exact
sequence of complexes
$$
0\to\underline\Gamma_Z\cI^\bullet\to\cI^\bullet
\to j_*j^*\cI^\bullet\to 0
$$
whence a natural exact sequence of $\cA$-modules :
\set\begin{equation}\label{eq_first-two-terms}
0\to\underline\Gamma_Z\cF\to\cF\to j_*j^*\cF\to
R^1\underline\Gamma_Z\cF\to 0
\qquad\text{for every $\cA$-module $\cF$}
\end{equation}
and natural isomorphisms :
\set\begin{equation}\label{eq_Gamma-as-a-j}
R^{q-1} j_*j^*\cF\isom R^q\underline\Gamma_Z\cF
\qquad\text{for all $q>1$}.
\end{equation}

\begin{lemma}\label{lem_flabby-Gamma_Z}
In the situation of \eqref{subsec_Gamma_Z}, let
$\cB\to\cA$ be any morphism of sheaves of rings on $X$.
The following holds :
\begin{enumerate}
\item
Every flabby $\cA$-module is $\underline\Gamma_Z$-acyclic.
\item
The functor $R\underline\Gamma_Z$ commutes with the
forgetful functor $\sD^+(\cA\Mod)\to\sD^+(\cB\Mod)$.
\item
Suppose that $X$ is locally coherent and quasi-separated,
and that $Z$ is a constructible closed subset of $X$. Then :
\begin{enumerate}
\item
Every qc-flabby $\cA$-module is $\underline\Gamma_Z$-acyclic.
\item
For every $i\in\N$ the functor $R^i\underline\Gamma_Z$
commutes with filtered colimits of $\cA$-modules.
\end{enumerate}
\end{enumerate}
\end{lemma}
\begin{proof} (i): First of all, if $\cF$ is flabby
on $X$, the natural morphism $\cF\to j_*j^*\cF$ is
obviously an epimorphism; together with
\eqref{eq_first-two-terms}, this implies that
$R^1\underline\Gamma_Z\cF=0$ when $\cF$ is flabby.
Next, \eqref{eq_Gamma-as-a-j}, together with the
fact that flabby $\cA$-modules are acyclic for direct
image functors such as $j_*$ (lemma \ref{lem_franziska}(ii)),
implies that $R^i\underline\Gamma_Z\cF=0$ also for $i>1$,
whence the contention.

Assertion (iii.a) is checked in the same way : first
one notices that, under the stated assumptions, the
morphism $\cF\to j_*j^*\cF$ is an epimorphism also
for all qc-flabby sheaves on $X$, which again gives
the vanishing of $R^1\underline\Gamma_Z\cF$; for
the higher derived functors one then appeals to
\eqref{eq_Gamma-as-a-j} and lemma \ref{lem_franziska}(iv).

(iii.b): Suppose first that $i>1$; in this case,
the assertion follows from \eqref{eq_Gamma-as-a-j}
and the following more general :

\begin{claim}\label{cl_gosh}
In the situation of (iii), the functor $R^ij_*$
commutes with filtered colimits, for every $i\in\N$.
\end{claim}
\begin{pfclaim} The claim can be checked at the stalks
at the points $x\in X$, hence we may assume that $X$
is quasi-compact, in which case the same holds for $U$.
Now, on the one hand, $R^qj_*j^*\cF$ is the sheaf
associated with the presheaf given by the rule :
$U'\mapsto H^q(U\cap U',\cF)$ for every open
subset $U'\subset X$, so $(R^qj_*j^*\cF)_x$ is
also the stalk at $x$ of this presheaf. On the other
hand, $x$ admits a fundamental system of open neighborhoods
consisting of quasi-compact open subsets, and if $U'$
is such a neighborhood of $x$, the subset $U'\cap U$
is also quasi-compact, and the functor
$\cF\mapsto H^q(U\cap U',\cF)$ commutes with filtered
colimits, by proposition \ref{prop_dir-im-and-colim}(ii).
The claim follows by combining these two observations.
\end{pfclaim}

In case $i\leq 1$, we argue similarly with
\eqref{eq_first-two-terms}, which reduces to
checking that both functors $\cF\mapsto\cF_x$ and
$\cF\mapsto(j_*j^*\cF)_x$ commute with filtered
colimits; this is trivial for the former functor,
and for the latter it is a special case of claim
\ref{cl_gosh}.

Lastly, (ii) is an immediate consequence of (i) and
remark \ref{rem_acyclic-crit}(iv), since any injective
$\cA$-module is a flabby $\cB$-module (lemma
\ref{lem_franziska}(v)).
\end{proof}

\sset\subsubsection{}\label{subsec_patch-trick}
In the situation of \eqref{subsec_Gamma_Z}, let $X'$
be the amalgamated sum in the following cocartesian diagram
of topological spaces :
$$
\xymatrix{ U \ar[r]^-j \ar[d]_j & X \ar[d]^{j'_1} \\
X \ar[r]^-{j'_2} & X'
}$$
and let $\pi:X'\to X$ be the unique continuous map
with $\pi\circ j'_i=\one_X$ for $i=1,2$. We notice :

\begin{proposition}\label{prop_patch-trick}
In the situation of \eqref{subsec_patch-trick}, let also
$\cF$ be any sheaf on $X$. Then we have :

{\em(i)}\ \
If $\cF$ is flabby, the same holds for $\pi^*\cF$.

{\em(ii)}\ \
If moreover $\cF$ is an $\cA$-module, the following holds :
\begin{enumerate}
\alphaenu
\item
We have a natural isomorphism in $\sD^+(\cA(X)\Mod)$ :
$$
R\Gamma\cF\oplus R\Gamma_{\!Z}\cF\isom R\Gamma(\pi^*\cF).
$$
\item
If $X$ is coherent and $Z$ is constructible in $X$,
then $R^i\Gamma_{\!Z}\cF=0$ for every $i>\dim X$.
\item
If $X$ is spectral and $Z$ is constructible in $X$, define
$d_\cF$ as in theorem {\em\ref{th_tohoku-vanish}(ii)}. Then
$$
R^i\Gamma_{\!Z}\cF=0
\qquad
\text{for every $i>d_\cF$}.
$$
\end{enumerate}
\end{proposition}
\begin{proof}(i): Set $V_i:=j'_i(X)$ for $i=1,2$; then $X'$ admits
the open covering $X'=V_1\cup V_2$, and the restriction of
$\pi^*\cF$ is flabby on $V_1$ and $V_2$. Then the assertion
follows from lemma \ref{lem_flabby-is-local}.

(ii.a): For $i=1,2$ and every $\cA$-module $\cG$, we have a
natural isomorphism of $\cA(X)$-modules
$$
\psi_{\cG,i}:\Gamma(X,\cF)\isom\Gamma(V_i,\cG).
$$
Moreover, for every $(\sigma,\tau)\in\Gamma\cG\oplus\Gamma_{\!Z}\cG$
there exists a unique section $\phi_\cG(\sigma,\tau)\in\Gamma(\pi^*\cG)$
whose restriction to $V_1$ agrees with $\psi_{\cG,1}(\sigma)$ and
whose restriction to $V_2$ agrees with $\psi_{\cG_2}(\sigma+\tau)$
(details left to the reader). The rule
$(\sigma,\tau)\mapsto\phi_\cG(\sigma,\tau)$ yields a natural
isomorphism
$$
\phi_\cG:\Gamma\cG\oplus\Gamma_{\!Z}\cG\isom\Gamma(\pi^*\cG).
$$
Now, let $\cF\to\cG^\bullet$ be a resolution of $\cF$ by a complex
of injective $\cA$-modules; by lemma \ref{lem_franziska}(v), the
complex $\pi^*\cG^\bullet$ is a resolution of $\pi^*\cF$ by flabby
$\pi^*\cA$-modules. Combining with lemma \ref{lem_flabby-Gamma_Z}(i)
and remark \ref{rem_acyclic-crit}(iv), we deduce that
$\phi_{\cG^\bullet}:\Gamma\cG^\bullet\oplus\Gamma_{\!Z}\cG^\bullet
\isom\Gamma(\pi^*\cG^\bullet)$ is the sought natural isomorphism.

Lastly, if $X$ is coherent (resp. spectral) and $Z$ is
constructible, the open subsets $V_1$ and $V_2$ are coherent
(resp. spectral) and retro-compact in $X'$; by lemma
\ref{lem_charact-coh-sob-sp}(iv), in this case $X'$ is
coherent (resp. $X'$ is spectral, and it is easily seen
that $d_{\pi^*\cF}=d_\cF$). Then (ii.b) and (ii.c) follow
directly from (ii.a) and theorem \ref{th_tohoku-vanish}.
\end{proof}

\sset\subsubsection{}\label{subsec_complement}
In the situation of \eqref{subsec_lambdasmus},
suppose that the maps $g$ and $g_\lambda$, for every
$\lambda\in\Ob(\Lambda)$ are open immersions, and set
$$
Z:=X\!\setminus\!Y
\qquad
Z_\lambda:=X_\lambda\!\setminus\!Y_\lambda
\qquad
\text{for every $\lambda\in\Ob(\Lambda)$}.
$$
Consider also a system
$(\cG_\lambda~|~\lambda\in\Ob(\Lambda))$, where
$\cG_\lambda$ is a $\Z_{X_\lambda}$-module for every
$\lambda\in\Lambda$, with transition maps
$$
\phi^{-1}_u\cG_\lambda\to\cG_\mu
\qquad
\text{for every morphism $u:\mu\to\lambda$ in $\Lambda$}
$$
and set
$$
\cG:=\colim_{\lambda\in\Ob(\Lambda)}\phi^*_\lambda\cG_\lambda.
$$
Notice that
$$
\phi^{-1}_uZ_\lambda\subset Z_\mu
\qquad\text{and}\qquad
\phi^{-1}_\lambda Z_\lambda\subset Z
$$
for every $\lambda\in\Ob(\Lambda)$ and every morphism
$u:\mu\to\lambda$ in $\Lambda$. There follows a system
of natural morphisms
$$
\tau_u:\phi^*_u\underline\Gamma_{Z_\lambda}\cG_\lambda\to
\underline\Gamma_{Z_\mu}\cG_\mu
\qquad\text{and}\qquad
\tau_\lambda:\phi^*_\lambda\underline\Gamma_{Z_\lambda}\cG_\lambda\to
\underline\Gamma_Z\cG
$$
amounting to a co-cone $\tau_\bullet$ with vertex
$\underline\Gamma_Z\cG$.

\begin{proposition}\label{prop_better-late}
In the situation of \eqref{subsec_complement}, the
co-cone $\tau_\bullet$ induces a natural isomorphism
$$
\colim_{\lambda\in\Ob(\Lambda)}
\phi^*_\lambda\,R^i\underline\Gamma_{Z_\lambda}\cG_\lambda
\isom R^i\underline\Gamma_Z\cG
\qquad
\text{for every $i\in\N$}.
$$
\end{proposition}
\begin{proof} From proposition \ref{prop_dir-im-and-colim}(ii)
we deduce a natural isomorphism
$$
\begin{aligned}
\colim_{\lambda\in\Ob(\Lambda)}
\phi^*_\lambda\,\underline\Gamma_{Z_\lambda}\cG_\lambda=\, &
\colim_{\lambda\in\Ob(\Lambda)}\phi^*_\lambda\,
\Ker\,(\cG_\lambda\to g_{\lambda*}g^*_\lambda\cG_\lambda) \\
\isom\, & \Ker\,(\cG\to g_*g^*\cG)=\underline\Gamma_Z\cG.
\end{aligned}
$$
Now, let us choose a compatible system of injective
resolutions $\cG_\lambda\to\cI^\bullet_\lambda$, so that
the transition maps for the system
$(\cG_\lambda~|~\lambda\in\Ob(\Lambda))$ extend to a
system of transition morphisms of complexes
$(\cI^\bullet_\lambda~|~\lambda\in\Ob(\Lambda))$.
On the one hand, claim \ref{cl_old-flame} says that
$$
\cI^\bullet:=
\colim_{\lambda\in\Ob(\Lambda)}\phi^*_\lambda\cI^\bullet_\lambda
$$
is a qc-flabby resolution of $\cG$. On the other hand,
the foregoing yields a natural isomorphism
$$
\colim_{\lambda\in\Ob(\Lambda)}
\phi^*_\lambda R^i\underline\Gamma_{Z_\lambda}\cG_\lambda
\isom H^i\underline\Gamma_Z\cI^\bullet.
$$
Then the assertion follows from lemma
\ref{lem_flabby-Gamma_Z}(iii.a).
\end{proof}

\begin{corollary}\label{cor_better-late}
Let $(X,\cA)$ be a ringed topological space with $X$
locally spectral, $Z\subset X$ a constructible closed
subset, $x\in Z$ any point, and denote by $j_x:X(x)\to X$
the inclusion map (notation of definition
{\em\ref{def_special}(iii)}). We have a natural
isomorphism of functors
$$
j_x^*R\underline\Gamma_Z\isom
R\underline\Gamma_{Z(x)}
\qquad
\text{in $\sD^+(j_x^*\cA\Mod)$}.
$$
\end{corollary}
\begin{proof} We may assume that $X$ is spectral,
in which case $X(x)$ is the limit of the cofiltered
system of all quasi-compact open neighborhoods of
$x$ in $X$, and if $U$ is any such neighborhood,
the inclusion map $U\to X$ is quasi-compact and
quasi-separated. The sought morphism $\phi^\bullet$
is then deduced from a co-cone constructed as in
\eqref{subsec_complement}, and it suffices to check
that $H^i\phi^\bullet$ is an isomorphism on the
underlying $\Z_{X(x)}$-modules for every $i\in\Z$,
so the assertion is reduced to proposition
\ref{prop_better-late}.
\end{proof}

\sset\subsubsection{}\label{subsec_alternate-Cech}
In the situation of \eqref{subsec_Gamma_Z}, notice the
natural isomorphism of $\cA$-modules :
$$
\cHom_\Z(i_*\Z_Z,\cF)\isom\underline\Gamma_Z\cF
\qquad
\text{for every $\cA$-module $\cF$}
$$
whence a natural isomorphism of $\cA(X)$-modules :
$$
\Hom_\Z(i_*\Z_Z,\cF)\isom\Gamma_{\!Z}\cF
\qquad
\text{for every $\cA$-module $\cF$}.
$$
Now, let $\fU:=(U_i~|~i\in I)$ be any open covering of $U$,
indexed by a set $I$, and fix any total ordering on $I$. By
proposition \ref{prop_Cech-resolve}(iv), we get natural
isomorphisms in $\sD^+(\cA\Mod)$ and $\sD^+(\cA(X)\Mod)$
$$
R\cHom^\bullet_\Z(R^\alt_\bullet(\fU)[1],K^\bullet)\isom
R\underline\Gamma_ZK^\bullet
\qquad
R\Hom^\bullet_\Z(R^\alt_\bullet(\fU)[1],K^\bullet)\isom
R\Gamma_{\!Z}K^\bullet
$$
for every bounded below complex $K^\bullet$ of $\cA$-modules.
These isomorphisms allow us to compute $R\Gamma_{\!Z}K^\bullet$
as a {\em \v{C}ech cohomology functor}, as follows.
Clearly, the rule : $\cF\mapsto C^\bullet_\mathrm{alt}(\fU,\cF)$
defines a functor $\cA\Mod\to\sC(\cA(X)\Mod)$ (notation of
definition \ref{def_Cech-complex}(ii)).
Now, pick any resolution $K^\bullet\isom\cI^\bullet$ by a bounded
below complex of injective $\cA$-modules; the foregoing yields a
natural isomorphism :
$$
R\Gamma_{\!Z}K^\bullet\isom
\Tot\,C^\bullet_\mathrm{alt}(\fU,\cI^\bullet)[-1]
\qquad
\text{in $\sD^+(\cA(X)\Mod)$}.
$$

\sset\subsubsection{}
We define a bifunctor
$$
\cHom^\bullet_Z:=\underline\Gamma_Z\circ\cHom^\bullet_\cA:
\sC(\cA\Mod)\times\sC(\cA\Mod)^o\to\sC(\cA\Mod).
$$
The bifunctor $\cHom^\bullet_Z$ admits a right derived functor :
$$
R\cHom^\bullet_Z:\sD^+(\cA\Mod)\times\sD(\cA\Mod)^o\to
\sD(\cA\Mod).
$$
The construction can be outlined as follows. First, for
a fixed complex $K_\bullet$ of $\cA$-modules, one can
consider the right derived functor of the functor
$\cG\mapsto\cHom^\bullet_Z(K_\bullet,\cG)$, which is denoted
$R\cHom_Z^\bullet(K_\bullet,-):\sD^+(\cA\Mod)\to\sD(\cA\Mod)$.
Next, one verifies that every quasi-isomorphism
$K_\bullet\to K'_\bullet$ induces an isomorphism
of functors :
$R\cHom^\bullet_Z(K'_\bullet,-)\isom R\cHom^\bullet_Z(K_\bullet,-)$,
hence the natural transformation
$K_\bullet\mapsto R\cHom^\bullet_Z(K_\bullet,-)$ factors through
$\sD(\cA\Mod)$, and this is the sought bifunctor
$R\cHom^\bullet_Z$. In case $Z=X$, one recovers the functor
$R\cHom^\bullet_{\cA}$ of \eqref{subsec_def-Hom-cplx}.

\begin{lemma}\label{lem_adj-Gamma_Z}
In the situation of \eqref{subsec_Gamma_Z}, the following holds :
\begin{enumerate}
\item
The functor $\underline\Gamma_Z$ is right adjoint to $i_*i^*$.
More precisely, for any two $\cA$-modules $\cF$ and $\cG$
we have natural isomorphisms :
\begin{enumerate}
\item
$\cHom_\cA(\cF,\underline\Gamma_Z\cG)\isom
\cHom_\cA(i_*i^*\cF,\cG)$.
\item
$R\cHom^\bullet_\cA(\cF,R\underline\Gamma_Z\cG)\isom
R\cHom^\bullet_\cA(i_*i^*\cF,\cG)$.
\end{enumerate}
\item
If $\cF$ is an injective (resp. flabby) $\cA$-module on $X$,
then the same holds for $\underline\Gamma_Z\cF$.
\item
There are natural isomorphisms of bifunctors :
$$
R\cHom^\bullet_Z\isom R\underline\Gamma_Z\circ R\cHom^\bullet_\cA
\isom R\cHom^\bullet_\cA(-,R\underline\Gamma_Z-)
$$
\item
If\/ $W\subset X$ is any other closed subset, there is a
natural isomorphism of functors :
$$
R\underline\Gamma_{W\cap Z}\isom
R\underline\Gamma_W\circ R\underline\Gamma_Z.
$$
\item
If $f:(Y,\cB)\to(X,\cA)$ is any morphism of
ringed spaces, there is a natural isomorphism of functors :
$$
Rf_*\circ R\underline\Gamma_{f^{-1}Z}\isom
R\underline\Gamma_Z\circ Rf_*:\sD^+(\cB\Mod)\to\sD^+(\cA\Mod).
$$
\item
There is a natural isomorphism of bifunctors on
$\sD^+(\cA\Mod)^o\times\sD^+(\cA\Mod)$ :
$$
R\cHom^\bullet_Z(R\underline\Gamma_Z-,-)\isom
R\cHom^\bullet_\cA(R\underline\Gamma_Z-,-).
$$
\end{enumerate}
\end{lemma}
\begin{proof} (i): To establish the isomorphism (a), one uses
the short exact sequence
$$
0\to j_!j^*\cF\to\cF\to i_*i^*\cF\to 0
$$
to show that any map $\cF\to\underline\Gamma_Z\cG$ factors
through $\cH:=i_*i^*\cF$. Conversely, since $j_*j^*\cH=0$,
it is clear that every map $\cH\to\cG$ must factor through
$\underline\Gamma_Z\cG$. The isomorphism (b) is derived easily
from (a) and (ii).

(iii): According to \cite[Th.10.8.2]{We}, the first isomorphism in
(iii) is deduced from lemmata \ref{lem_flabby-Gamma_Z}(i) and
\ref{lem_triv-dual}(ii). The second isomorphism in (iii) follows
from \cite[Th.10.8.2]{We}, and assertion (ii) concerning
injective sheaves.

(ii): Suppose first that $\cF$ is injective; let $\cG_1\to\cG_2$
a monomorphism of $\cA$-modules and $\phi:\cG_1\to\underline\Gamma_Z\cF$
any $\cA$-linear map. By (i) we deduce a map $i_*i^*\cG_1\to\cF$,
which extends to a map $i_*i^*\cG_2\to\cF$, by injectivity of $\cF$.
Finally, (i) again yields a map $\cG_2\to\underline\Gamma_Z\cF$
extending $\phi$.

Next, if $\cF$ is flabby, let $V\subset X$ be any open subset,
and $s\in\underline\Gamma_Z\cF(V)$; since $s_{|V\setminus Z}=0$,
we can extend $s$ to a section
$s'\in\underline\Gamma_Z\cF(V\cup(X\!\setminus\! Z))$.
Since $\cF$ is flabby, $s'$ extends to a section $s''\in\cF(X)$;
however, by construction $s''$ vanishes on the complement
of $Z$.

(iv): Clearly
$\underline\Gamma_{W\cap Z}=\underline\Gamma_W\circ\underline\Gamma_Z$,
so the claim follows easily from (ii) and \cite[Th.10.8.2]{We}.

(v): Since flabby $\cB$-modules are acyclic for $f_*$ and
$\underline\Gamma_{f^{-1}Z}$ (lemmata \ref{lem_franziska}(ii)
and \ref{lem_flabby-Gamma_Z}(i)), one can apply (ii) and
\cite[Th.10.8.2]{We} to identify both functors with
$R(\underline\Gamma_Z\circ f_*)$.

(vi): Let $I^\bullet$ and $J^\bullet$ be two bounded below
complexes of injective $\cA$-modules. Clearly any morphism
$\underline\Gamma_ZI^\bullet\to J^\bullet$ vanishes
on $U$, whence the assertion.
\end{proof}

\begin{definition} Let $X$ be a topological space, and
$K^\bullet\in\sD^+(\Z_X\Mod)$.

(i)\ \
The {\em depth of $K^\bullet$ along $Z$} is
$$
\depth_Z K^\bullet:=
\sup\{n\in\Z~|~R^i\underline\Gamma_ZK^\bullet=0\text{ for all $i<n$}\}
\in\Z\cup\{+\infty\}.
$$
From lemma \ref{lem_adj-Gamma_Z}(ii,iv) and the Grothendieck
spectral sequence (\cite[Th.5.8.3]{Mat}), we see that
\set\begin{equation}\label{eq_spec-depth}
\depth_W K^\bullet\geq\depth_Z K^\bullet
\qquad \text{whenever $W\subset Z$}.
\end{equation}

(ii)\ \
Let $\Phi:=(Z_\lambda~|~\lambda\in\Lambda)$ be a family of closed
subsets of $X$, and suppose that $\Phi$ is cofiltered by inclusion.
Due to \eqref{eq_spec-depth} it is reasonable to define the
{\em depth of $K^\bullet$ along $\Phi$} as
$$
\depth_\Phi K^\bullet:=
\sup\,\{\depth_{Z_\lambda} K^\bullet~|~\lambda\in\Lambda\}
\in\Z\cup\{+\infty\}.
$$
\end{definition}

\sset\subsubsection{}\label{subsec_case-of-schs}
Let now $X$ be a scheme and $i:Z\to X$ a closed immersion of
schemes, such that $Z$ is constructible in $X$. Then the open
immersion $j:X\!\setminus\! Z\to X$ is quasi-compact and
separated, hence the functors $R^qj_*$ preserve the subcategories
of quasi-coherent modules, for every $q\in\N$
(\cite[Ch.III, Prop.1.4.10]{EGAIII}). In light of
\eqref{eq_Gamma-as-a-j} we deduce that $R^q\underline\Gamma_Z$
restricts to a functor :
$$
R^q\underline\Gamma_Z:\cO_{\!X}\Mod_\qcoh\to\cO_{\!X}\Mod_\qcoh
\qquad\text{for every $q\in\N$}.
$$

\begin{lemma}\label{lem_without-cohereur}
In the situation of \eqref{subsec_case-of-schs}, the following
holds :
\begin{enumerate}
\item
The functor $R\underline\Gamma_Z$ restricts
to a triangulated functor :
$$
R\underline\Gamma_Z:
\sD^+(\cO_{\!X}\Mod)_\qcoh\to\sD^+(\cO_{\!X}\Mod)_\qcoh
$$
(notation of \eqref{subsec_der-qcoh}).
\item
Let $f:Y\to X$ be an affine morphism of schemes, $K^\bullet$
any object of\/ $\sD^+(\cO_Y\Mod)_\qcoh$. Then we have the identity :
$$
\depth_{f^{-1}Z}K^\bullet=\depth_Z Rf_*K^\bullet.
$$
\item
Let $f:Y\to X$ be a flat morphism of schemes, $K^\bullet$
any object of\/ $\sD^+(\cO_{\!X}\Mod)_\qcoh$. Then the
natural map
$$
f^*R\underline\Gamma_ZK^\bullet\to
R\underline\Gamma_{f^{-1}Z}f^*K^\bullet
$$
is an isomorphism in $\sD^+(\cO_Y\Mod)_\qcoh$. Especially,
we have
$$
\depth_{f^{-1}Z}f^*K^\bullet\geq\depth_ZK^\bullet
$$
and if $f$ is faithfully flat, the inequality is actually an
equality.
\end{enumerate}
\end{lemma}
\begin{proof} (i): Indeed, if $K^\bullet$ is an object of
$\sD^+(\cO_{\!X}\Mod)_\qcoh$, we have a spectral sequence :
$$
E_2^{pq}:=R^p\underline\Gamma_Z H^qK^\bullet\Rightarrow
R^{p+q}\underline\Gamma_Z K^\bullet.
$$
(This spectral sequence is obtained from a Cartan-Eilenberg
resolution of $K^\bullet$ : see {\em e.g.} \cite[\S5.7]{We}.)
Hence $R^\bullet\underline\Gamma_Z K^\bullet$ admits
a finite filtration whose subquotients are quasi-coherent;
then the lemma follows from \cite[Ch.III, Prop.1.4.17]{EGAIII}.

(ii): To start with, a spectral sequence argument as in the
foregoing shows that $Rf_*$ restricts to a functor
$\sD^+(\cO_Y\Mod)_\qcoh\to\sD^+(\cO_{\!X}\Mod)_\qcoh$, and
moreover we have natural identifications :
$R^if_*K^\bullet\isom f_*H^iK^\bullet$
for every $i\in\Z$ and every object $K^\bullet$ of
$\sD^+(\cO_Y\Mod)_\qcoh$. In view of lemma \ref{lem_adj-Gamma_Z}(v),
we derive natural isomorphisms :
$R^i\underline\Gamma_Z Rf_*K^\bullet\isom
Rf_*R^i\underline\Gamma_{f^{-1}Z}K^\bullet$ for every
object $K^\bullet$ of $\sD^+(\cO_Y\Mod)_\qcoh$.
The assertion follows easily.

(iii): A spectral sequence argument as in the proof of
(i) allows us to assume that $K^\bullet=\cF[0]$ for some
quasi-coherent $\cO_{\!X}$-module $\cF$. In this case, let
$j_X:X\!\setminus\!Z\to X$ and $j_Y:Y\!\setminus\!f^{-1}Z\to Y$
be the open immersions; considering the exact sequences 
$f^*\eqref{eq_first-two-terms}$ and $f^*\eqref{eq_Gamma-as-a-j}$,
we reduce to showing that the natural map
$f^*Rj_{X*}j_X^*\cF\to Rj_{Y*}j_Y^*f^*\cF$ is
an isomorphism in $\sD^+(\cO_Y\Mod)_\qcoh$. The latter
assertion follows from corollary \ref{cor_base-change-where}.
\end{proof}

\begin{proposition}\label{prop_affine-Cech-resolve}
Let $X$ be an affine scheme, $Z\subset X$ any closed subset,
$\cF^\bullet$ a bounded below complex of quasi-coherent
$\cO_{\!X}$-modules, and $\fU:=(U_t~|~t\in I)$ a family
of affine open subsets of $X$ such that
$\bigcup_{t\in I}U_t=X\!\setminus\!Z$. There exists a natural
isomorphism
$$
R\Gamma_{\!Z}\cF^\bullet\isom
\Tot\,C_\mathrm{alt}^\bullet(\fU,\cF^\bullet)[-1]
\qquad
\text{in $\sD^+(\cO_{\!X}(X)\Mod)$}
$$
(where $C_\mathrm{alt}^\bullet$ denotes the augmented
alternating \v{C}ech complex).
\end{proposition}
\begin{proof} Pick any resolution $\cF^\bullet\isom\cI^\bullet$
by a bounded below complex of injective $\cO_{\!X}$-modules; by
virtue of \eqref{subsec_alternate-Cech} it suffices to show that
the natural map of double complexes
$$
C_\mathrm{alt}^\bullet(\fU,\cF^\bullet)\to
C_\mathrm{alt}^\bullet(\fU,\cI^\bullet)
$$
induces a quasi-isomorphism on total complexes. To this aim, we
can argue as in the proof of theorem \ref{th_Cech-resolve}(ii) :
the details shall be left to the reader (the condition that $X$
is affine, is needed here, since we are dealing with the
augmented \v{C}ech complex, whereas in the proof of theorem
\ref{th_Cech-resolve}(ii) only the ordinary version intervenes,
so for the proof of the latter it is only required that the
intersection of any finite system of members of $\fU$ is affine).
\end{proof}

\sset\subsubsection{}\label{subsec_local-depth}
Let $X$ be a scheme; for every $x\in X$, we consider the
cofiltered family $\Phi(x)$ of all non-empty constructible
closed subsets of $X(x):=\Spec\,\cO_{\!X,x}$; clearly
$\bigcap_{Z\in\Phi(x)}Z=\{x\}$. Let $K^\bullet$ be any object
of $\sD^+(\Z_X\Mod)$; we are interested in the quantities :
$$
\delta(x,K^\bullet):=\depth_{\{x\}} K^\bullet(x)
\qquad\text{and}\qquad
\delta'(x,K^\bullet):=\depth_{\Phi(x)} K^\bullet(x)
$$
(notation of definition \ref{def_strict-loc}(iii)). Especially,
we wish to compute the depth of a complex $K^\bullet$ along a
closed subset $Z\subset X$, in terms of the local invariants
$\delta(x,K^\bullet)$ or $\delta'(x,K^\bullet)$, evaluated at the
points $x\in Z$. This shall be achieved by theorem
\ref{th_local-depth}. To begin with, we remark :

\begin{lemma}\label{lem_ineq-deltas}
With the notation of \eqref{subsec_local-depth}, we have :
\begin{enumerate}
\item
$\delta(x,K^\bullet)\geq\delta'(x,K^\bullet)$\ \ for every
$x\in X$ and every $K^\bullet\in\Ob(\sD^+(\Z_X\Mod))$.
\item
If the topological space $|X|$ underlying $X$ is locally noetherian,
the inequality of\/ {\em (i)} is actually an equality.
\item
If $\cF$ is a flat quasi-coherent $\cO_{\!X}$-module, then
$\delta'(x,\cF)\geq\delta'(x,\cO_{\!X})$.
\end{enumerate}
\end{lemma}
\begin{proof} (i) follows easily from \eqref{eq_spec-depth}.
Likewise, (ii) follows likewise from \eqref{eq_spec-depth}
and the fact that if $|X|$ is locally noetherian, $\{x\}$ is
the smallest element of the family $\Phi(x)$.

(iii): According to \cite[Ch.I, Th.1.2]{La}, the
$\cO_{\!X(x)}$-module $\cF(x)$ is the colimit of a
filtered system $(\cL_\lambda~|~\lambda\in\Lambda)$ of free
$\cO_{\!X(x)}$-modules of finite rank. In view of lemma
\ref{lem_flabby-Gamma_Z}(iii.b), it is easily seen that
$$
\depth_Z \cF\geq
\inf_{\lambda\in\Lambda}\depth_Z \cL_\lambda=\depth_Z\cO_{\!X}
$$
for every closed onstructible subset $Z\subset X$. The
assertion follows.
\end{proof}

\begin{theorem}\label{th_local-depth}
With the notation of \eqref{subsec_local-depth}, let
$Z\subset X$ be any closed constructible subset, $K^\bullet$
any object of\/ $\sD^+(\Z_X\Mod)$. Then 
$$
\depth_ZK^\bullet=\inf\,\{\delta(x,K^\bullet)~|~x\in Z\}=
\inf\,\{\delta'(x,K^\bullet)~|~x\in Z\}.
$$
\end{theorem}
\begin{proof} Let $d:=\depth_ZK^\bullet$; in light of
corollary \ref{cor_better-late}, it is clear that
$\delta(x,K^\bullet)\geq d$ for all $x\in Z$, hence,
in order to prove the first identity it suffices
to show :
\begin{claim} Suppose $d<+\infty$. Then there exists $x\in Z$
such that $R^d\Gamma_{\!\{x\}}K^\bullet(x)\neq 0$.
\end{claim}
\begin{pfclaim} By definition of $d$, we can find an open
subset $V\subset X$ and a non-zero section
$s\in\Gamma(V,R^d\underline\Gamma_ZK^\bullet)$. The support
of $s$ is a closed subset $S\subset Z$.  Let $x$ be a maximal
point of $S$. From lemma \ref{lem_adj-Gamma_Z}(iv) we deduce a
spectral sequence
$$
E_2^{pq}:=R^p\underline\Gamma_{\{x\}}
          R^q\underline\Gamma_{Z\cap X(x)}K^\bullet(x)
\Rightarrow R^{p+q}\underline\Gamma_{\{x\}}K^\bullet(x)
$$
and corollary \ref{cor_better-late} implies that
$E_2^{pq}=0$ whenever $q<d$, therefore :
$$
R^d\underline\Gamma_{\{x\}}K^\bullet(x)\simeq
R^0\underline\Gamma_{\{x\}}
R^d\underline\Gamma_{Z\cap X(x)}K^\bullet(x).
$$
However, the image of $s$ in
$R^d\underline\Gamma_{Z\cap X(x)}K^\bullet(x)$ is supported
precisely at $x$, as required.
\end{pfclaim}

Finally, since $Z$ is constructible, $Z\cap X(x)\in\Phi(x)$ for
every $x\in X$, hence $\delta'(x,K^\bullet)\geq d$. Then the
second identity follows from the first and lemma
\ref{lem_ineq-deltas}(i).
\end{proof}

\begin{corollary}\label{cor_local-depth}
Let $Z\subset X$ be a closed and constructible subset, $\cF$
a $\Z_X$-module, and $j:X\!\setminus\! Z\to X$ the induced
open immersion. Then the natural map $\cF\to j_*j^*\cF$ is an
isomorphism if and only if\/ $\delta(x,\cF)\geq 2$, if and only
if\/ $\delta'(x,\cF)\geq 2$ for every $x\in Z$.
\qed\end{corollary}

\sset\subsubsection{}\label{subsec_compute-depth}
Let now $X$ be an affine scheme, say $X:=\Spec\,A$ for some ring
$A$, and $Z\subset X$ a constructible closed subset. In this case,
we wish to show that the depth of a complex of quasi-coherent
$\cO_{\!X}$-modules along $Z$ can be computed in terms of $\Ext$
functors on the category of $A$-modules. This is the content
of the following :

\begin{proposition}\label{prop_compute-depth}
In the situation of \eqref{subsec_compute-depth}, let $N$ be
any finitely presented $A$-module such that $\Supp\,N=Z$, and
$M^\bullet$ any bounded below complex of $A$-modules. Then :
$$
\depth_Z M^{\bullet\sim}=
\sup\,\{n\in\Z~|~\Ext^j_A(N,M^\bullet)=0\text{ for all $j<n$}\}.
$$
\end{proposition}
\begin{proof} ($M^{\bullet\sim}$ is the complex of quasi-coherent
$\cO_{\!X}$-modules determined by $M^\bullet$, as in
\eqref{sec_various-O-mod}).

\begin{claim}\label{cl_loc-to-glob-Ext}
There is a natural isomorphism in $\sD^+(A\Mod)$ :
$$
R\Hom^\bullet_A(N,M^\bullet)\isom
R\Hom^\bullet_{\cO_{\!X}}(N^\sim,R\underline\Gamma_ZM^{\bullet\sim}).
$$
\end{claim}
\begin{pfclaim} Let $i:Z\to X$ be the closed immersion.
The assumption on $N$ implies that $N^\sim=i_*i^*N^\sim$;
hence lemma \ref{lem_adj-Gamma_Z}(i.b) yields a natural isomorphism
\set\begin{equation}\label{eq_from-supp-N}
R\Hom^\bullet_{\cO_{\!X}}(N^\sim,M^{\bullet\sim})\isom
R\Hom^\bullet_{\cO_{\!X}}(N^\sim,R\underline\Gamma_Z M^{\bullet\sim}).
\end{equation}
The claim follows by combining corollary \ref{cor_Ext-loc=glob}(i)
and \eqref{eq_from-supp-N}.
\end{pfclaim}

From claim \ref{cl_loc-to-glob-Ext} we see already that
$\Ext^j_A(N,M^\bullet)=0$ whenever $j<\depth_ZM^{\bullet\sim}$.
Suppose that $\depth_ZM^{\bullet\sim}=p<+\infty$; in this case,
claim \ref{cl_loc-to-glob-Ext} gives an isomorphism :
$$
\Ext^p_A(N,M^\bullet)\simeq
\Hom_{\cO_{\!X}}(N^\sim,R^p\underline\Gamma_ZM^{\bullet\sim})\simeq
\Hom_A(N,R^p\Gamma_{\!\!Z}M^{\bullet\sim})
$$
where the last isomorphism holds by lemma \ref{lem_without-cohereur}.
To conclude the proof, we have to exhibit a non-zero
map from $N$ to $Q:=R^p\Gamma_{\!\!Z}M^{\bullet\sim}$. Let $F_0(N)\subset A$
denote the Fitting ideal of $N$; this is a finitely generated
ideal, whose zero locus coincides with the support of $N$. More precisely,
$\Ann_A(N)^r\subset F_0(N)\subset\Ann_A(N)$ for all sufficiently
large integers $r>0$. Let now $x\in Q$ be any non-zero element,
and $f_1,\dots,f_k$ a finite system of generators for $F_0(N)$;
since $x$ vanishes on $X\!\setminus\! Z$, the image of $x$ in
$Q_{f_i}$ is zero for every $i\leq k$, {\em i.e.} there exists
$n_i\geq 0$ such that $f_i^{n_i}x=0$ in $Q$. It follows easily
that $F_0(N)^n\subset\Ann_A(x)$ for a sufficiently large integer
$n$, and therefore $x$ defines a map $\phi:A/F_0(N)^n\to Q$ by
the rule : $a\mapsto ax$. According to \cite[Lemma 3.2.21]{Ga-Ra},
we can find a finite filtration
$0=J_m\subset\cdots\subset J_1\subset J_0:=A/F_0(N)^n$
such that each $J_i/J_{i+1}$ is quotient of a direct sum
of copies of $N$. Let $s\leq m$ be the smallest integer such
that $J_s\subset\Ker\,\phi$. By restriction, $\phi$ induces
a non-zero map $J_{s-1}/J_s\to Q$, whence a non-zero map
$N^{(S)}\to Q$, for some set $S$, and finally a non-zero map
$N\to Q$, as required.
\end{proof}

\begin{remark} Notice that the existence of a finitely presented
$A$-module $N$ with $\Supp\,N=Z$ is equivalent to the
constructibility of $Z$.
\end{remark}

\sset\subsubsection{}\label{subsec_depth_A}
In the situation of \eqref{subsec_compute-depth}, let
$I\subset A$ be a finitely generated ideal such that
$V(I)=Z$, and $M^\bullet$ a bounded below complex of
$A$-modules. Then it is customary to set :
$$
\depth_IM^\bullet:=\depth_ZM^{\bullet\sim}
$$
and the depth along $Z$ of $M^{\bullet\sim}$ is also
called the {\em $I$-depth\/} of $M^\bullet$. Moreover,
if $A$ is a local ring with maximal ideal $\fm_A$,
we shall often use the standard notation :
$$
\depth_AM^\bullet:=\depth_{\{\fm_A\}}M^{\bullet\sim}
$$
and this invariant is often called briefly the {\em depth\/}
of $M^\bullet$. With this notation, theorem \ref{th_local-depth}
can be rephrased as the identity :
\set\begin{equation}\label{eq_for-I-depth}
\depth_IM^\bullet=\inf_{\fp\in V(I)}\depth_{A_\fp}M^\bullet_\fp.
\end{equation}

\sset\subsubsection{}\label{subsec_alt-Cech}
Suppose now that $\bff$ is a finite system of generators
of the ideal $I\subset A$, and set $Z:=\Spec_,A/I$. 
In this generality, the Koszul complex $\bK^\bullet(\bff)$
of remark \ref{rem_koszul-alg}(ii) is not a resolution of
the $A$-module $A/I$, and consequently $H^\bullet(\bff,M)$
(for an $A$-module $M$) does not necessarily agree with
$\Ext^\bullet_A(A/I,M)$. Nevertheless, the following holds.

\begin{proposition}\label{prop_depth-Kosz}
{\em(i)}\ \
In the situation of \eqref{subsec_alt-Cech}, let
$M^\bullet\in\Ob(\sD^+(A\Mod))$; then :
$$
\depth_IM^\bullet=d:=
\sup\,\{i\in\Z~|~H^j(\bff,M^\bullet)=0\text{ for all $j<i$}\}.
$$

{\em(ii)}\ \
If $d<+\infty$, there are natural $A$-linear isomorphisms :
$$
H^d(\bff,M^\bullet)\isom
\Hom_A(A/I,R^d\Gamma_{\!\!Z} M^{\bullet\sim})\isom
\Ext^d_A(A/I,M^\bullet).
$$

{\em(iii)}\ \
Moreover we have a natural isomorphism in $\sD(A\Mod)$ :
$$
\colim_{m\in\N}\bK^\bullet(\bff^m,M^\bullet)\isom
R\Gamma_{\!\!Z} M^{\bullet\sim}
$$
where the transition maps in the colimit are the maps
$\bphi^\bullet_{\bff^n}$ of \eqref{subsec_def-bphis}.
\end{proposition}
\begin{proof} (i): Let $B:=\Z[t_1,\dots,t_r]\to A$ be the
ring homomorphism defined by the rule : $t_i\mapsto f_i$ for
$i=1,\dots,r$; we denote $\psi:X\to Y:=\Spec\,B$ the induced
morphism, $J\subset B$ the ideal generated by the system
$\bt:=(t_i~|~i=1,\dots,r)$, and set $W:=V(J)\subset Y$.
From lemma \ref{lem_adj-Gamma_Z}(v) we deduce a natural isomorphism
in $\sD(\cO_Y\Mod)$ :
\set\begin{equation}\label{eq_depth-Z=W}
\psi_*R\underline\Gamma_ZM^{\bullet\sim}\isom
R\underline\Gamma_W\psi_*M^{\bullet\sim}.
\end{equation}
Hence we are reduced to showing :
\begin{claim}\label{cl_same-for-J} $d=\depth_J\psi_*M^\bullet$.
\end{claim}
\begin{pfclaim} Notice that $\bt$ is a regular sequence in $B$,
hence $\Ext_B^\bullet(B/J,\psi_*M^\bullet)\simeq
H^\bullet(\bt,\psi_*M^\bullet)$; by
proposition \ref{prop_compute-depth}, there follows the identity :
$$
\depth_J\psi_*M^\bullet=
\sup\,\{n\in\N~|~H^i(\bt,\psi_*M^\bullet)=0\text{ for all $i<n$}\}.
$$
Then the assertion follows from lemma \ref{lem_koszul-vanish}(iv).
\end{pfclaim}

(ii): From lemma \ref{lem_adj-Gamma_Z}(i.b) and corollary
\ref{cor_Ext-loc=glob}(i) we derive a natural isomorphism :
$$
R\Hom^\bullet_B(B/J,\psi_*M^\bullet)\isom
R\Hom_{\cO_Y}^\bullet((B/J)^\sim,R\underline\Gamma_W\psi_*M^{\bullet\sim}).
$$
However, due to \eqref{eq_depth-Z=W} and claim \ref{cl_same-for-J}
we may compute :
$$
R^d\Hom_{\cO_Y}^\bullet((B/J)^\sim,R\underline\Gamma_W\psi_*M^{\bullet\sim})
\isom\Hom_B(B/J,R^d\Gamma_{\!\!W}\psi_*M^{\bullet\sim})\isom
\Hom_A(A/I,R^d\Gamma_{\!\!Z}M^{\bullet\sim}).
$$
The first claimed isomorphism easily follows. Similarly, one applies
lemma \ref{lem_adj-Gamma_Z}(i.b) and corollary \ref{cor_Ext-loc=glob}(i)
to compute $\Ext^\bullet_A(A/I,M^\bullet)$, and deduces the second
isomorphism of (ii) using (i). Assertion (iii) follows immediately
from \eqref{eq_Koszul-Cech} and proposition
\ref{prop_affine-Cech-resolve}.
\end{proof}

\begin{corollary}\label{cor_depth-Kosz}
Let $A$ be a noetherian ring, $I\subset A$ an ideal,
$M^\bullet$ an object of $\sD^+(A\Mod)$, and set
$Z:=\Spec\,A/I$. We have a natural isomorphism
$$
\colim_{n\in\N}R^i\Hom^\bullet_A(A/I^n,M^\bullet)\isom
R^i\Gamma_{\!Z}M^{\bullet\sim}
\qquad
\text{for every $i\in\Z$}.
$$
\end{corollary}
\begin{proof} This is obtained by combining proposition
\ref{prop_depth-Kosz}(iii) and remark \ref{rem_Hartsho}.
\end{proof}

\begin{corollary}\label{cor_f-flat-invariance}
In the situation of \eqref{subsec_depth_A}, let
$B$ be a faithfully flat $A$-algebra. Then :
$$
\depth_{IB} B\otimes_AM^\bullet=\depth_IM^\bullet.
$$
\end{corollary}
\begin{proof} It is a special case of lemma
\ref{lem_without-cohereur}(iii). Alternatively, one remarks
that
$$
\bK^\bullet(\bg,M^\bullet\otimes_AB)\simeq
\bK^\bullet(\bg,M^\bullet)\otimes_AB
\qquad
\text{for every finite sequence of elements $\bg$ in $A$}
$$
which allows us to apply proposition \ref{prop_depth-Kosz}(iii).
\end{proof}

\begin{theorem}\label{th_depth-flat-basechange}
Let $A\to B$ be a local homomorphism of local rings, and suppose
that the maximal ideals $\fm_A$ and $\fm_B$ of $A$ and $B$ are
finitely generated. Let also $M$ be an $A$-module and $N$ a
$B$-module which is flat over $A$. Then we have :
$$
\depth_B(M\otimes_AN)=\depth_AM+\depth_{B/\fm_AB}(N/\fm_AN).
$$
\end{theorem}
\begin{proof} Let $\bff:=(f_1,\dots,f_r)$ and
$\bar\bg:=(\bar g_1,\dots,\bar g_s)$ be finite systems of generators
for $\fm_A$ and respectively $\bar\fm_B$, the maximal ideal of
$B/\fm_AB$. We choose an arbitrary lifting of $\bar\bg$ to a
finite system $\bg$ of elements of $\fm_B$; then $(\bff,\bg)$
is a system of generators for $\fm_B$. We have a natural
identification of complexes of $B$-modules :
$$
\bK_\bullet(\bff,\bg)\isom
\Tot_\bullet(\bK_\bullet(\bff)\otimes_A\bK_\bullet(\bg))
$$
whence natural isomorphisms :
$$
\bK^\bullet((\bff,\bg),M\otimes_AN)\isom
\Tot^\bullet(\bK^\bullet(\bff,M)\otimes_A\bK^\bullet(\bg,N)).
$$
A standard spectral sequence, associated with the filtration by rows,
converges to the cohomology of this total complex, and since $N$
is a flat $A$-module, $\bK^\bullet(\bg,N)$ is a complex of flat
$A$-modules, so that the $E_1$ term of this spectral sequence is
found to be :
$$
E^{pq}_1\simeq H^q(\bff,M)\otimes_A\bK^p(\bg,N)
$$
and the differentials $d_1^{pq}:E^{pq}_1\to E^{p+1,q}_1$ are
induced by the differentials of the complex $\bK^\bullet(\bg,N)$.
Set $\kappa_A:=A/\fm_A$; notice that $H^q(\bff,M)$ is a
$\kappa_A$-vector space, hence
$$
H^q(\bff,M)\otimes_AH^\bullet(\bg,N)\simeq
H^\bullet(H^q(\bff,M)\otimes_{\kappa_A}\bK^\bullet(\bar\bg,N/\fm_AN))
$$
and consequently :
$$
E^{pq}_2\simeq H^q(\bff,M)\otimes_{\kappa_A}H^p(\bar\bg,N/\fm_AN).
$$
Especially :
\set\begin{equation}\label{eq_def-d-d'}
\text{$E^{pq}_2=0$ if either $p<d:=\depth_AM$ or
$q<d':=\depth_{B/\fm_AB}(N/\fm_AN)$.}
\end{equation}
Hence $H^i((\bff,\bg),M\otimes_AN)=0$ for all $i<d+d'$
and moreover, if $d$ and $d'$ are finite,
$H^{d+d'}((\bff,\bg),M\otimes_AN)\simeq E^{dd'}_2\neq 0$.
Now the sought identity follows readily from proposition
\ref{prop_depth-Kosz}(i).
\end{proof}

\begin{remark}\label{rem_prod-of-Exts}
Let $d$ and $d'$ be as in \eqref{eq_def-d-d'}, and
suppose that $d$ and $d'$ are finite. By inspection
of the proof, we see that -- under the assumptions
of theorem \ref{th_depth-flat-basechange} -- a natural
isomorphism has been exhibited :
$$
\Ext^{d+d'}_A(B/\fm_B,M\otimes_AN)\simeq
\Ext^d_A(A/\fm_A,M)\otimes_A\Ext^{d'}_{B/\fm_AB}(B/\fm_B,N/\fm_AN).
$$
\end{remark}

\begin{definition} Let $(A,\fm)$ be a local ring, and
$M_\bullet\in\Ob(\sD^b(A\Mod))$. The {\em residual Tor dimension}
$\mathrm{rTor.dim}_A(M_\bullet)$ of $M_\bullet$ is
$\sup\{i\in\Z~|~H_i(M_\bullet\derotimes_AA/\fm)\neq 0\}\in
\Z\cup\{\pm\infty\}$.
\end{definition}

\begin{remark}
If $A$ is a noetherian local ring, and $M$ is any $A$-module
of finite type, then $\mathrm{rTor.dim}_A(M)$ equals the
projective dimension of $M$ : see \cite[\S19, Lemma 1]{Mat}.
\end{remark}

\begin{theorem} Let $A$ be a local noetherian ring, and
$M_\bullet\in\Ob(\sD^b(A\Mod))$. Then :

{\em(i)}\ \
$\depth_A(M_\bullet)=+\infty$ if and only if\/
$\mathrm{rTor.dim}_A(M_\bullet)=-\infty$.

{\em(ii)}\ \
If\/ $\mathrm{rTor.dim}_A(M_\bullet)\in\Z$, we have the
{\em Auslander-Buchsbaum identity} :
$$
\depth_A(M_\bullet)+\mathrm{rTor.dim}_A(M_\bullet)=\depth_AA.
$$
\end{theorem}
\begin{proof}(ii): Let $d:=\dim A$, and pick a
{\em system of parameters} $\bff$ of $A$, {\em i.e.} a
sequence $(f_1,\dots,f_d)$ of $d$ elements of $A$ that
generates a primary ideal. Using a Cartan-Eilenberg resolution
of the Koszul complex $\bK_\bullet(\bff)$ (see remark
\ref{rem_koszul-alg}(ii)), we get a spectral sequence (see
\cite[\S5.7]{We})
$$
E^2_{pq}:=H_p(M_\bullet\derotimes_AH_q(\bff,A))\Rightarrow
H_{p+q}(\bff,M_\bullet).
$$
Now, by proposition \ref{prop_depth-Kosz}(i) and lemma
\ref{lem_koszul-vanish}(v), we have :
$$
\begin{aligned}
\depth_AA&=d-e_A & & \text{where }
e_A:=\inf\{i\in\Z~|~H_j(\bff,A)=0\text{ for all $i>j$}\} \\
\depth_AM_\bullet&=d-e_{M_\bullet} & & \text{where }
e_{M_\bullet}:=
\inf\{i\in\Z~|~H_j(\bff,M_\bullet)=0\text{ for all $i>j$}\}.
\end{aligned}
$$
Thus, we have to show that if
$t_{M_\bullet}:=\mathrm{rTor.dim}_A(M_\bullet)\in\Z$, then
\set\begin{equation}\label{eq_corner}
e_{M_\bullet}=e_A+t_{M_\bullet}.
\end{equation}

\begin{claim}\label{cl_corner}
For every $A$-module $N\neq 0$ of finite length, we have :
$$
H_{t_{M_\bullet}}(M_\bullet\derotimes N)\neq 0
\qquad\text{and}\qquad
H_i(M_\bullet\derotimes N)=0
\qquad\text{for every $i>t_{M_\bullet}$}.
$$
\end{claim}
\begin{pfclaim} There exists a finite filtration
$0:=N_0\subset N_1\subset\cdots\subset N_{n-1}\subset N_n:=N$
such that $N_{i+1}/N_i$ is isomorphic to $A/\fm$ for
$i=0,\dots,n-1$. Hence, by a simple induction on the length
$n$ of the filtration, we are reduced to checking the following
assertion : let $\Sigma_\bullet:=(0\to N'\to N\to A/\fm\to 0)$
be a short exact sequence of $A$-modules of finite length. The
induced map
$H_i(M_\bullet\derotimes_AN')\to H_i(M_\bullet\derotimes_AN)$ is an
isomorphism for every $i\geq t_{M_\bullet}$. The latter follows
easily, y considering the long exact homology sequence induced
by the short exact sequence of complexes
$M_\bullet\derotimes_A\Sigma_\bullet$.
\end{pfclaim}

Now, since $A$ is noetherian and $\bff$ generates a primary ideal,
$H_q(\bff,A)$ is a module of finite length for every $q\in\N$ (lemma
\ref{lem_koszul-vanish}(i)). Claim \ref{cl_corner} implies that
$0\neq E^2_{t_{M_\bullet},e_A}$, and $E^2_{pq}=0$ whenever either
$p>t_{M_\bullet}$ or $q>e_A$, so that
$E^2_{t_{M_\bullet},e_A}=E^\infty_{t_{M_\bullet},e_A}$, whence \eqref{eq_corner}.

(i): If $t_{M_\bullet}\in\Z$, then \eqref{eq_corner} implies that
$e_{M_\bullet}\in\Z$, since $e_A\in\N$; then in this case
$\depth_AM_\bullet\in\Z$. If $t_{M_\bullet}=-\infty$, claim
\ref{cl_corner} says that $E^2_{pq}=0$ for every $p,q\in\Z$,
whence $H_i(\bff,M)=0$ for every $i\in\Z$, {\em i.e.}
$\depth_A(M_\bullet)=+\infty$. Lastly, by proposition
\ref{prop_depth-Kosz}(iii), this condition holds if and only
if $\colim_{m\in\N}\bK^\bullet(\bff^m,M^\bullet)=0$ in
$\sD^b(A\Mod)$. Set
$$
K_\bullet:=\colim_{m\in\N}\bK_\bullet(\bff^m)
\qquad
T_{m,\bullet}:=t^{\leq-1}\bK_\bullet(\bff^m)
\quad\text{for every $m\in\N$}\qquad
T_\bullet:=\colim_{m\in\N}T_{m,\bullet}.
$$
Then, in turn the latter means that
$M^\bullet\otimes_AK_\bullet=0$ in $\sD^b(A\Mod)$ (lemma
\ref{lem_koszul-vanish}(v)). Set $k:=A/\fm$; the short
exact of complexes $0\to A[0]\to K_\bullet\to T_\bullet\to 0$
induces a distinguished triangle in $\sD^-(A\Mod)$ :
$$
M_\bullet\derotimes_Ak[0]\to
M_\bullet\derotimes_AK_\bullet\derotimes_Ak[0]\to
M_\bullet\derotimes_AT_\bullet\derotimes_Ak[0]\to
M_\bullet\derotimes_Ak[1].
$$
Notice that the transition maps $T_{m,\bullet}\to T_{m+1,\bullet}$
induce the zero maps $T_{m,\bullet}\otimes_Ak\to T_{m+1,\bullet}\otimes_Ak$
for every $m\in\N$; hence $T_\bullet\derotimes_Ak[0]=0$ in
$\sD(A\Mod)$, and therefore
$M_\bullet\derotimes_AT_\bullet\derotimes_Ak[0]=0$ as well.
Summing up, we conclude that if $\depth_AM_\bullet=+\infty$,
then $M_\bullet\derotimes_Ak[0]=0$ in $\sD(A\Mod)$, {\em i.e.}
$t_{M^\bullet}=-\infty$.
\end{proof}

\sset\subsubsection{}
Let $f:X\to S$ be any morphism of schemes. Recall that an
$\cO_{\!X}$-module $\cF$ is said to be
{\em $f$-flat at a point $x\in X$\/} if $\cF_x$ is a flat
$\cO_{S,f(x)}$-module. We say that $\cF$ is 
{\em $f$-flat over a point $s\in S$\/} if $\cF$ is $f$-flat
at all points of $f^{-1}(s)$. Finally, we say that $\cF$
is {\em $f$-flat\/} if $\cF$ is $f$-flat at all the points
of $X$ (\cite[Ch.0, \S6.7.1]{EGAI}).

\begin{corollary}\label{cor_depth-flat-basechange}
Let $f:X\to S$ be a finitely presented morphism of schemes,
$\cF$ a $f$-flat quasi-coherent $\cO_{\!X}$-module of finite
presentation, $\cG$ any quasi-coherent $\cO_{\!S}$-module, $x\in X$
any point, $i:f^{-1}f(x)\to X$ the natural morphism. Then :
$$
\delta'(x,\cF\otimes_{\cO_{\!X}}f^*\cG)=
\delta'(x,i^*\cF)+\delta'(f(x),\cG).
$$
\end{corollary}
\begin{proof} Set $s:=f(x)$, and denote by
$$
j_s:S(s):=\Spec\,\cO_{\!S,s}\to S\qquad \text{and} \qquad
f_s:X(s):=X\times_SS(s)\to S(s)
$$
the induced morphisms. By inspecting the definitions, one
checks easily that
$$
\delta'(x,\cF\otimes_{\cO_{\!X}}f^*\cG)=
\delta'(x,\cF_{|X(s)}\otimes_{\cO_{\!X(s)}}f_s^*\cG(s))
\qquad \text{and} \qquad
\delta'(s,\cG)=\delta'(s,\cG(s))
$$
(notation of definition \ref{def_strict-loc}(iii)).
Thus, we may replace $S$ by $S(s)$ and $X$ by $X(s)$, and
assume from start that $S$ is a local scheme and $s$ its closed
point. Clearly we may also assume that $X$ is affine and
finitely presented over $S$. Then we can write $S$ as the
limit of a cofiltered family $(S_\lambda~|~\lambda\in\Lambda)$
of local noetherian schemes, with local transition maps, and
we may assume that $f:X\to S$ is a limit of a cofiltered family
$(f_\lambda:X_{\!\lambda}\to S_\lambda~|~\lambda\in\Lambda)$
of morphisms of finite type, such that
$X_\mu=X_\lambda\times_{S_\lambda}S_\mu$ whenever $\mu\geq\lambda$.
Likewise, we may descend $\cF$ to a family
$(\cF_\lambda~|~\lambda\in\Lambda)$ of finitely presented
quasi-coherent $\cO_{\!X_{\!\lambda}}$-modules, such that
$\cF_\mu=\phi^*_{\mu\lambda}\cF_\lambda$ for every $\mu\geq\lambda$
(where $\phi_{\mu\lambda}:X_\mu\to X_{\!\lambda}$ is the natural
morphism). Furthermore, up to replacing $\Lambda$ by some
cofinal subset, we may assume that $\cF_\lambda$ is
$f_\lambda$-flat for every $\lambda\in\Lambda$
(\cite[Ch.IV, Cor.11.2.6.1(ii)]{EGAIV-3}).
For every $\lambda\in\Lambda$, denote by
$\psi_\lambda:S\to S_\lambda$ the natural morphism, so that
$s_\lambda:=\psi_\lambda(s)$ is the closed point of $S_\lambda$;
for every constructible closed subset
$Z_\lambda\subset S_\lambda$, the set $Z:=\psi^{-1}_\lambda Z_\lambda$
is constructible and closed in $S$, and according to lemma
\ref{lem_without-cohereur}(ii) we have
$\depth_Z\cG=\depth_{Z_\lambda}\psi_{\lambda*}\cG$.
On the other hand, every closed constructible subset
$Z\subset S$ is of the form $\psi^{-1}_\lambda Z_\lambda$ for
some $\lambda\in\Lambda$ and some $Z_\lambda$ closed constructible
in $S_\lambda$ (\cite[Ch.IV, Cor.8.3.11]{EGAIV-3}), so that :
\set\begin{equation}\label{eq_about-G}
\delta'(s,\cG)=\sup\{\delta'(s_\lambda,\psi_{\lambda*}\cG)~|~
\lambda\in\Lambda\}.
\end{equation}
Similarly, for every $\lambda\in\Lambda$ let $x_{\!\lambda}$
be the image of $x$ in $X_{\!\lambda}$ and denote by
$i_\lambda:f^{-1}_\lambda(s_\lambda)\to X_{\!\lambda}$ the natural
morphism. Notice that
$f^{-1}_\lambda(s_\lambda)$ is a scheme of finite type over
$\Spec\,\kappa(s_\lambda)$, especially it is noetherian,
hence $\delta'(x_\lambda,i^*_\lambda\cF_\lambda)=
\delta(x_\lambda,i^*_\lambda\cF_\lambda)$, by lemma
\ref{lem_ineq-deltas}(ii); by the same token we have as
well the identity : $\delta'(x,i^*\cF)=\delta(x,i^*\cF)$.
Denote by $g_\lambda:f^{-1}(s)\to f^{-1}_\lambda(s_\lambda)$
and $g_{\lambda,x}:\Spec\,\cO_{\!f^{-1}(s),x}\to
\Spec\,\cO_{\!f^{-1}_\lambda(s_\lambda),x_\lambda}$ the induced morphisms.

\begin{claim}\label{cl_yorik}
There exists $\lambda\in\Lambda$ such that
$g_{\mu,x}^{-1}(x_\mu)=\{x\}$ for every $\mu\geq\lambda$.
\end{claim}
\begin{pfclaim} Let $Z$ be the topological closure of $\{x\}$
in $f^{-1}(s)$; we may find $\lambda\in\Lambda$ and a closed
subset $Z_\lambda\subset f^{-1}_\lambda(s_\lambda)$ such that
$Z=g_\lambda^{-1}Z_\lambda$ (\cite[Ch.IV, Cor.8.3.11]{EGAIV-3}).
Since $Z$ is irreducible and $g_\lambda$ is surjective, it
is easily seen that $Z_\lambda$ is irreducible.
Moreover, the topological closure $Z'_\lambda$ of $\{x_\lambda\}$
in $f_\lambda^{-1}(s_\lambda)$ lies in $Z_\lambda$, and
$g_\lambda^{-1}Z'_\lambda$ contains $x$; thus, $Z_\lambda=Z'_\lambda$.
It follows that $g_{\lambda,x}^{-1}(x_\lambda)=\{x\}$, whence the
assertion.
\end{pfclaim}

Since $g_{\lambda,x}$ is faithfully flat, from claim \ref{cl_yorik}
and corollary \ref{cor_f-flat-invariance} we see that, after
replacing $\Lambda$ by a cofinal subset, we may assume :
\set\begin{equation}\label{eq_about-F}
\delta'(x_{\!\lambda},i^*_\lambda\cF_\lambda)=\delta'(x,i^*\cF)
\qquad\text{for every $\lambda\in\Lambda$}.
\end{equation}

\begin{claim}\label{cl_general-over-K}
Let $K$ be a field, $X$ and $Y$ two $K$-schemes, and set
$Z:=X\times_{\Spec\,K}Y$. Let also $p:Z\to X$ and $q:Z\to Y$
be the projections, $\cM$ any $\cO_{\!X}$-module, and $\cN$
any $\cO_Y$-module. Then we have :
$$
\Supp\,(p^*\cM\otimes_{\cO_{\!Z}}q^*\cN)=
p^{-1}(\Supp\,\cM)\cap q^{-1}(\Supp\,\cN).
$$
\end{claim}
\begin{pfclaim} To ease notation, set
$\cP:=p^*\cM\otimes_{\cO_{\!Z}}q^*\cN$. Let $z\in Z$ be any
point with $x:=p(z)\in\Supp\,\cM$ and $y:=q(z)\in\Supp\,\cN$;
it suffices to show that $\cP_z\neq 0$. Now, $\cO_{Z,z}$ is
a localization of $\cO_{\!X,x}\otimes_K\cO_{Y,y}$, at a prime
ideal $\fp$, and $\cP_z=(\cM_x\otimes_K\cN_y)_\fp$. We come
down therefore to the following situation. Let $A$ and $B$
be two local $K$-algebras, $\fm_A\subset A$,
$\fm_B\subset B$ the respective maximal ideals,
$M\neq 0$ an $A$-module, $N\neq 0$ a $B$-modules,
and $\fp\subset A\otimes_KB$ a prime ideal such that
\set\begin{equation}\label{eq_all-right}
\fm_A\otimes_KB+A\otimes_K\fm_B\subset\fp.
\end{equation}
We have to check that $(M\otimes_KN)_\fp\neq 0$. To this
aim, we may reduce easily to the case where $M$ and $N$
are cyclic modules, in which case there exist a surjective
$A$-linear map $M\to\kappa(A):=A/\fm_A$ and a surjective
$B$-linear map $N\to\kappa(B):=B/\fm_B$; we are then further
reduced to showing that
$(\kappa(A)\otimes_K\kappa(B))_\fp\neq 0$, which is ensured
by \eqref{eq_all-right}.
\end{pfclaim}

Next, notice that the local scheme $X(x)$ is the limit of the
cofiltered system of local schemes
$(X_{\!\lambda}(x_{\!\lambda})~|~\lambda\in\Lambda)$.
Let $Z_\lambda\subset X_{\!\lambda}(x_{\!\lambda})$ be a closed
constructible subset and set $Z:=\phi_\lambda^{-1}Z_\lambda$,
where $\phi_\lambda:X(x)\to X_{\!\lambda}(x_{\!\lambda})$ is the
natural morphism.
Let $Y_\lambda$ be the fibre product in the cartesian diagram
of schemes :
$$
\xymatrix{ Y_\lambda \ar[r]^-{q_\lambda} \ar[d]_{p_\lambda} &
           X_{\!\lambda}(x_{\!\lambda}) \ar[d] \\
           S \ar[r]^-{\psi_\lambda} & S_\lambda.
}$$
The morphism $\phi_\lambda$ factors through a unique morphism
of $S$-schemes $\bar\phi_\lambda:X(x)\to Y_\lambda$,
and if $y_\lambda\in Y_\lambda$ is the image of the closed
point $x$ of $X(x)$, then $\bar\phi_\lambda$ induces a natural
identification $X(x)\isom Y_\lambda(y_\lambda)$. To ease notation,
let us set :
$$
\cH:=\cF\otimes_{\cO_{\!X}}f^*\cG
\qquad\text{and}\qquad
\cH_\lambda:=
\cF_\lambda\otimes_{\cO_{\!X_{\!\lambda}}}f_\lambda^*\psi_{\lambda*}\cG.
$$
In view of corollary \ref{cor_better-late}, there follows
a natural isomorphism :
\set\begin{equation}\label{eq_here-we-go}
R\underline\Gamma_Z\cH(x)\isom
\bar\phi{}^*_\lambda R\underline\Gamma_{q^{-1}_\lambda(Z_\lambda)}
(q^*_\lambda\cF_\lambda(x_\lambda)\otimes_{\cO_{Y_\lambda}}p_\lambda^*\cG).
\end{equation}
On the other hand, applying the projection formula (see remark
\ref{rem_dir-im-and-colim}) we get :
\set\begin{equation}\label{eq_proj-form}
q_{\lambda*}(q^*_\lambda\cF_\lambda(x_\lambda)
\otimes_{\cO_{Y_\lambda}}p_\lambda^*\cG)\simeq
\cH_\lambda(x_\lambda).
\end{equation}
Combining \eqref{eq_here-we-go}, \eqref{eq_proj-form} and
lemma \ref{lem_adj-Gamma_Z}(v) we deduce :
$$
R\underline\Gamma_Z\cH(x)\simeq
\phi^*_\lambda R\underline\Gamma_{Z_\lambda}
\cH_\lambda(x_{\!\lambda}).
$$
Notice that the foregoing argument applies also in case $S$ is
replaced by some $S_\mu$ for some $\mu\geq\lambda$ (and
consequently $X$ is replaced by $X_\mu$); we arrive therefore
at the inequality :
\set\begin{equation}\label{eq_deltas}
\delta'(x_\lambda,\cH_\lambda)\leq\delta'(x_\mu,\cH_\mu)
\leq\delta'(x,\cH)
\qquad\text{whenever $\mu\geq\lambda$}.
\end{equation}

\begin{claim}\label{cl_hard-deltas} $\delta'(x,\cH)=
\sup\,\{\delta'(x_\lambda,\cH_\lambda)~|~\lambda\in\Lambda\}$.
\end{claim}
\begin{pfclaim} Due to \eqref{eq_deltas} we may assume --
again after replacing $\Lambda$ by a cofinal subset --
that $\delta'(x_\lambda,\cH_\lambda)$ is a constant $d\in\N$
independent of $\lambda$. Especially,
$R^d\underline\Gamma_{Z_\lambda}\cH_\lambda(x_{\!\lambda})\neq 0$
for every $\lambda\in\Lambda$ and every closed constructible
subset $Z_\lambda\subset X_{\!\lambda}(x_{\!\lambda})$. We have to
show that $R^d\underline\Gamma_Z\cH(x)\neq 0$ for arbitrarily
small constructible closed subsets $Z\subset X(x)$. However, given
such $Z$, there exist $\lambda\in\Lambda$ and $Z_\lambda$ closed
constructible in $X_{\!\lambda}(x_{\!\lambda})$ such that
$Z=\phi^{-1}(Z_\lambda)$ (\cite[Ch.IV, Cor.8.3.11]{EGAIV-3}).
Say that:
$$
S=\Spec\,A \qquad S_\lambda=\Spec\,A_\lambda \qquad
X_{\!\lambda}(x_{\!\lambda})=\Spec\,B_\lambda.
$$
Hence $Y_\lambda\simeq\Spec\,A\otimes_{A_\lambda}B_\lambda$.
Let also $F_\lambda$ (resp. $G$) be a $B_\lambda$-module
(resp. $A$-module) such that $F^\sim_\lambda\simeq\cF_\lambda$
(resp. $G^\sim\simeq\cG$). Then
$\cH_{|Y_\lambda}\simeq(F_\lambda\otimes_{A_\lambda}G)^\sim$.
Let $\fm_\lambda\subset B_\lambda$ and
$\fn_\lambda\subset A_\lambda$ be the maximal ideals.
Up to replacing $Z$ by a smaller subset, we may assume
that $Z_\lambda=V(\fm_\lambda)$. To ease notation, set :
$$
E_1:=\Ext^a_{A_\lambda}(\kappa(s_\lambda),G)
\qquad
E_2:=\Ext^b_{B_\lambda/\fn_\lambda B_\lambda}
(\kappa(x_\lambda),F_\lambda/\fn_\lambda F_\lambda)
$$
with $a=\delta'(s_\lambda,\psi_{\lambda*}\cG)$ and $b=d-a$.
Proposition \ref{prop_depth-Kosz}(i),(ii) and remark
\ref{rem_prod-of-Exts} say that
$$
0\neq E_3:=\Ext^d_{B_\lambda}(\kappa(x_\lambda),F_\lambda\otimes_{A_\lambda}G)
\simeq E_1\otimes_{\kappa(s_\lambda)}E_2.
$$
To conclude, it suffices to show that $y_\lambda\in\Supp\,E_3$.
However, claim \ref{cl_general-over-K} says that
$$
\Supp\,E_3=p_\lambda^{-1}(\Supp\,E_1)\cap q_\lambda^{-1}(\Supp\,E_2).
$$
Now, $p_\lambda(y_\lambda)=s$, and clearly $s\in\Supp\,E_1$,
since $G$ is an $A$-module and $s$ is the closed point of $S$;
likewise, $q_\lambda(y_\lambda)=x_\lambda$ is the closed point of
$X_\lambda(x_\lambda)$, therefore it is in the support of $E_2$.
\end{pfclaim}

We can now conclude the proof : indeed, from theorem
\ref{th_depth-flat-basechange} we derive :
$$
\delta'(x_\lambda,\cH_\lambda)=
\delta'(x,i^*_\lambda\cF_\lambda)+\delta'(f(x),\psi_{\lambda*}\cG)
\qquad\text{for every $\lambda\in\Lambda$}
$$
and then the corollary follows from \eqref{eq_about-G},
\eqref{eq_about-F} and claim \ref{cl_hard-deltas}.
\end{proof}

\sset\subsubsection{}\label{subsec_fibrewise-depth}
Consider now a finitely presented morphism of schemes $f:X\to Y$,
and a sheaf of ideals $\cI\subset\cO_{\!X}$ of finite type. Let $x\in X$
be any point, and set $y:=f(x)$; notice that $f^{-1}(y)$ is an algebraic
$\kappa(y)$-scheme, so that the local ring $\cO_{\!f^{-1}(y),x}$ is
noetherian. We let $\cI(y):=i^{-1}_y\cI\cdot\cO_{\!f^{-1}(y)}$, which
is a coherent sheaf of ideals of the fibre $f^{-1}(y)$.
The {\em fibrewise $\cI$-depth of $X$ over $Y$ at the point $x$\/}
is the integer:
$$
\depth_{\cI,f}(x):=\depth_{\cI(y)_x}\cO_{\!f^{-1}(y),x}.
$$
In case $\cI=\cO_{\!X}$, the fibrewise $\cI$-depth shall be called
the {\em fibrewise depth\/} at the point $x$, and shall be denoted
by $\depth_f(x)$. In view of \eqref{eq_for-I-depth} we have the
identity:
\set\begin{equation}\label{eq_for-fibrewise-I-dp}
\depth_{\cI,f}(x)=\inf\{\depth_f(z)~|~z\in V(\cI(y))\text{ and }
x\in\overline{\{z\}}\}.
\end{equation}

\sset\subsubsection{}
The fibrewise $\cI$-depth can be computed locally as follows.
For a given $x\in X$, set $y:=f(x)$ and let $U\subset X$ be
any affine open neighborhood of $x$ such that $\cI_{|U}$ is generated
by a finite family $\bff:=(f_i)_{1\leq i\leq r}$, where $f_i\in\cI(U)$
for every $i\leq r$. Then proposition \ref{prop_depth-Kosz}(i)
implies that :
\set\begin{equation}\label{eq_local-Koszul}
\depth_{\cI,f}(x)=\sup\{n\in\N~|~
H^i(\bK^\bullet(\bff,\cO_{\!X}(U))_x\otimes_{\cO_{Y,y}}\kappa(y))=0
\text{ for all } i<n\}.
\end{equation}

\begin{proposition}
In the situation of \eqref{subsec_fibrewise-depth}, for every integer
$d\in\N$ we define the subset:
$$
L_\cI(d):=\{x\in X~|~\depth_{\cI,f}(x)\geq d\}.
$$
Then the following holds:
\begin{enumerate}
\item
$L_\cI(d)$ is a constructible subset of $X$ (see
definition {\em\ref{def_constructible}(v)}).
\item
Let $Y'\to Y$ be any morphism of schemes, set $X':=X\times_YY'$,
and let $g:X'\to X$, $f':X'\to Y'$ be the induced morphisms. Then:
$$
L_{g^*\cI}(d)=g^{-1}L_\cI(d)\qquad\text{for every $d\in\N$.}
$$
\end{enumerate}
\end{proposition}
\begin{proof} (ii): The identity can be checked on the fibres,
hence we may assume that $Y=\Spec\,k$ and $Y'=\Spec\,k'$ for some
fields $k\subset k'$. Let $x'\in X'$ be any point and $x\in X$
its image; since the map $\cO_{\!X,x}\to\cO_{\!X',x'}$ is faithfully
flat, corollary \ref{cor_f-flat-invariance} implies:
$$
\depth_{g^*\cI,f'}(x')=\depth_{\cI,f}(x)
$$
whence (ii).

(i): We may assume that $Y=\Spec\,A$,
$X=\Spec\,B$, $\cI=I^\sim$, where $A$ is a ring, $B$
is a finitely presented $A$-algebra, and $I\subset B$
is a finitely generated ideal, for which we choose a
finite system of generators $\bff:=(f_i)_{1\leq i\leq r}$.
We may also assume that $A$ is reduced.
Suppose first that $A$ is noetherian; then, for every $i\in\N$
and $y\in Y$ let us set:
$$
N_\cI(i)=\bigcup_{y\in Y}
\Supp\,H^i(\bK^\bullet(\bff,B)\otimes_A\kappa(y)).
$$
Taking into account \eqref{eq_for-I-depth} and
\eqref{eq_local-Koszul}, we deduce:
$$
L_\cI(d)=
\bigcap^{d-1}_{i=0}(X\!\setminus\! N_\cI(i)).
$$
It suffices therefore to show that each subset $N_\cI(i)$ is
constructible. For every $j\in\N$ we set:
$$
\bZ^j:=\Ker(d^j:\bK^j(\bff,B)\to\bK^{j+1}(\bff,B))\qquad
\bB^j:=\Img(d^{j-1}:\bK^{j-1}(\bff,B)\to\bK^j(\bff,B)).
$$
Using \cite[Ch.IV, Cor. 8.9.5]{EGAIV-3}, one deduces easily
that there exists an affine open subscheme $U\subset Y$,
say $U=\Spec\,A'$ for some flat $A$-algebra $A'$, such that
$A'\otimes_A\bZ^\bullet$, $A'\otimes_A\bB^\bullet$ and
$A'\otimes_A H^\bullet(\bK^\bullet(\bff,B))$ are flat
$A'$-modules. By noetherian induction, we can then replace
$Y$ by $U$ and $X$ by $X\times_YU$, and assume from start
that $\bZ^\bullet$, $\bB^\bullet$ and
$H^\bullet(\bK^\bullet(\bff,B))$ are flat $A$-modules.
In such case, taking homology of the complex
$\bK^\bullet(\bff,B)$ commutes with any base change; therefore:
$$
N_\cI(i)=\bigcup_{y\in Y}
\Supp\, H^i(\bff,B)\otimes_A\kappa(y)=
\Supp\, H^i(\bff,B)
$$
whence the claim, since the support of a $B$-module of finite
type is closed in $X$.

Finally, for a general ring $A$, we can find a noetherian
subalgebra $A'\subset A$, an $A'$-algebra $B'$ of finite
type and a finitely generated ideal $I'\subset B'$ such that
$B=A\otimes_{A'}B$ and $I=I'B$. Let $\cI'$ be the sheaf of
ideals on $X':=\Spec\,B'$ determined by $I'$; by the foregoing,
$L_{\cI'}(d)$ is a constructible subset of $X'$. Thus, the
assertion follows from (ii), and the fact that a morphism
of schemes is continuous for the constructible topology
(\cite[Ch.IV, Prop.1.8.2]{EGAIV}).
\end{proof}

\subsection{Depth and associated primes}

\begin{definition}\label{def_Ass}
Let $X$ be a scheme, $\cF$ a quasi-coherent $\cO_{\!X}$-module.
\begin{enumerate}
\item
Let $x\in X$ be any point; we say that $x$ is {\em associated with\/}
$\cF$ if there exists $f\in\cF_x$ such that the radical of the
annihilator of $f$ in $\cO_{\!X,x}$ is the maximal ideal of
$\cO_{\!X,x}$. If $x$ is a point associated with $\cF$ and $x$
is not a maximal point of $\Supp\,\cF$, we say that $x$ is an
{\em imbedded point\/} for $\cF$. We shall denote :
$$
\Ass\,\cF:=\{x\in X~|~\text{$x$ is associated with $\cF$}\}.
$$
\item
We say that the $\cO_{\!X}$-module $\cF$
{\em satisfies condition\/} $S_1$ if every associated point
of $\cF$ is a maximal point of $X$.
\item
Likewise, if $X$ is affine, say $X=\Spec\,A$, and $M$ is
any $A$-module, we denote by $\Ass_A\,M\subset X$ (or just
$\Ass\,M$, if the notation is not ambiguous) the set of
prime ideals associated with the $\cO_{\!X}$-module $M^\sim$
arising from $M$. An associated (resp. imbedded) point of
$M^\sim$ is also called an {\em associated\/} (resp. {\em imbedded})
{\em prime ideal\/} of $M$.

We say that $M$ {\em satisfies condition\/}
$S_1$ if the same holds for the $\cO_{\!X}$-module $M^\sim$.
\item
Let $x\in X$ be a point, $\cG$ a quasi-coherent $\cO_{\!X}$-submodule
of $\cF$. We say that $\cG$ is a {\em $x$-primary
submodule\/} of $\cF$ if $\Ass\,\cF/\cG=\{x\}$. We say that
$\cG$ is a {\em primary submodule \/} of $\cF$ if there exists a
point $x\in X$ such that $\cG$ is $x$-primary.
\item
We say that a submodule $\cG$ of the $\cO_{\!X}$-module $\cF$
admits a {\em primary decomposition\/} if there exist
primary submodules $\cG_1,\dots,\cG_n\subset\cF$ such that
$\cG=\cG_1\cap\cdots\cap\cG_n$.
\end{enumerate}
\end{definition}

\begin{remark}
(i)\ \ 
Our definition of associated point is borrowed
from \cite[Partie I, D{\'e}f.3.2.1]{Gr-Ra}, and it corresponds
to Bourbaki's notion of {\em weakly associated prime}
(\cite[Ch.IV, \S1, Exerc.17]{BouAC}), also called
{\em weak Bourbaki prime} in other works. Bourbaki's definition
of associated prime, in \cite[Ch.IV, \S1, no.1]{BouAC}
(which is the same as in \cite[Exp.III, D{\'e}f.1.1]{SGA2}),
agrees with ours for modules over noetherian rings (see lemma
\ref{lem_ass-are-mins}(ii)), but in general the two notions
diverge.

(ii)\ \ 
For a noetherian ring $A$ and a finitely generated $A$-module
$M$, our condition $S_1$ is the same as in \cite{Har3}.
It also agrees with that of \cite[Ch.IV, D\'ef.5.7.2]{EGAIV-2},
in case $\Supp\,M=\Spec\,A$.
\end{remark}

\begin{lemma}\label{lem_ass-are-mins}
Let $A$ be a ring, and $M$ any $A$-module. The following holds :
\begin{enumerate}
\item
$\Ass\,M$ is the set of all $\fp\in\Spec\,A$ with the following
property. There exists $m\in M$, such that $\fp$ is a maximal
point of the closed subset $\Supp(m)$ ({\em i.e.} $\fp$ is the
preimage of a minimal prime ideal of the ring $A/\Ann_A(m)$).
\item
If $\fp\in\Ass\,M$, and $\fp$ is finitely generated, there exists
$m\in M$ with $\fp=\Ann_A(m)$.
\end{enumerate}
\end{lemma}
\begin{proof}(i): Indeed, suppose that $m\in M\!\setminus\!\{0\}$
is any element, and $\fp$ is a maximal point of $V(\Ann_A(m))$;
denote by $m_\fp\in M_\fp:=M\otimes_AA_\fp$ the image of $m$. Then :
\set\begin{equation}\label{eq_klar}
\Ann_{A_\fp}(m_\fp)=\Ann_A(m)\otimes_AA_\fp.
\end{equation}
Hence $A_\fp/\Ann_{A_\fp}(m_\fp)=(A/\Ann_A(m))\otimes_AA_\fp$,
and the latter is by assumption a local ring of Krull dimension
zero; it follows easily that the radical of $\Ann_{A_\fp}(m_\fp)$
is $\fp A_\fp$, {\em i.e.} $\fp\in\Ass_AM$, as stated.

Conversely, suppose that $\fp\in\Ass_AM$; then there exists
$m_\fp\in M_\fp$ such that the radical of $\Ann_{A_\fp}(m_\fp)$
equals $\fp A_\fp$. We may assume that $m_\fp$ is the image
of some $m\in M$; then \eqref{eq_klar} implies that
$V(\Ann_A(m))\cap\Spec\,A_\fp=V(\Ann_{A_\fp}(m_\fp))=\{\fp\}$,
so $\fp$ is a maximal point of $V(\Ann_A(m))$.

(ii): Suppose $\fp\in\Ass\,M$ is finitely generated; by definition,
there exists $y\in M_\fp$ such that the radical of $I:=\Ann_{A_\fp}(y)$
equals $\fp A_\fp$.
Since $\fp$ is finitely generated, some power of $\fp A_\fp$ is
contained in $I$. Let $t\in\N$ be the smallest integer such
that $\fp^tA_\fp\subset I$, and pick any non-zero element $x$
in $\fp^{t-1}A_\fp y$; then $\Ann_{A_\fp}(x)=\fp A_\fp$.
After clearing some denominators, we may assume that $x\in M$.
Let $a_1,\dots,a_r$ be a finite set of generators for $\fp$;
we may then find $t_1,\dots,t_r\in A\!\setminus\!\fp$ such
that $t_ia_ix=0$ in $M$, for every $i\leq r$. Set
$x':=t_1\cdot{\dots}\cdot t_rx$; then $\fp\subset\Ann_A(x')$;
however $\Ann_{A_\fp}(x')=\fp A_\fp$, therefore $\fp=\Ann_A(x')$,
whence the contention.
\end{proof}

\begin{lemma}\label{lem_Ass-Supp}
Let $X$ be a scheme, $\cF$ a quasi-coherent $\cO_{\!X}$-module. Then :
\begin{enumerate}
\item
$\Ass\,\cF\subset\Supp\,\cF$.
\item
$\Ass\,\cF_{|U}=U\cap\Ass\,\cF$\/ for every open subset\/ $U\subset X$.
\item
$\cF=0$\/ if and only if\/ $\Ass\,\cF=\emptyset$.
\end{enumerate}
\end{lemma}
\begin{proof} (i) and (ii) are obvious. To show (iii), we may assume
that $\cF\neq 0$, and we have to prove that $\Ass\,\cF$ is
not empty. Let $x\in X$ be any point such that $\cF_x\neq 0$,
and pick a non-zero section $f\in\cF_x$; set
$I:=\Ann_{\cO_{\!X,x}}(f)$, let $\bar\fq$ be a minimal prime
ideal of the quotient ring $\cO_{\!X,x}/I$, and
$\fq\subset\cO_{\!X,x}$ the preimage of $\bar\fq$. Then
$\fq$ corresponds to a generization $y$ of the point $x$
in $X$; let $f_y\in\cF_y$ be the image of $f$. Then
$\Ann_{\cO_{\!X,y}}(f_y)=I\cdot\cO_{Y,y}$, whose radical
is the maximal ideal $\fq\cdot\cO_{\!X,y}$, {\em i.e.}
$y\in\Ass\,\cF$.
\end{proof}

\begin{proposition}\label{prop_misc-Ass}
Let $X$ be a scheme, $\cF$ a quasi-coherent $\cO_{\!X}$-module.
We have :
\begin{enumerate}
\item
$\Ass\,\cF=\{x\in X~|~\delta(x,\cF)=0\}$.
\item
If\ \ $0\to\cF'\to\cF\to\cF''$\ \  is an exact sequence of
quasi-coherent $\cO_{\!X}$-modules, then :
$$
\Ass\,\cF'\subset\Ass\,\cF\subset\Ass\,\cF'\cup\Ass\,\cF''.
$$
\item
If $x\in X$ is any point, and $\cG_1,\dots,\cG_n$
(for some $n>0$) are $x$-primary submodules of $\cF$,
then $\cG_1\cap\cdots\cap\cG_n$ is also $x$-primary.
\item
Let $f:X\to Y$ be a finite morphism of schemes. Then
\begin{enumerate}
\item
The natural map :
$$
\bigoplus_{x\in f^{-1}(y)}\Gamma_{\!\{x\}}\cF_{|X(x)}\to
\Gamma_{\!\{y\}}f_*\cF_{|Y(y)}
$$
is a bijection for all $y\in Y$.
\item
$\Ass\,f_*\cF=f(\Ass\,\cF)$.
\end{enumerate}
\item
Suppose that $\cF$ is the union of a filtered family
$(\cF_\lambda~|~\lambda\in\Lambda)$ of quasi-coherent
$\cO_{\!X}$-submodules. Then :
$$
\bigcup_{\lambda\in\Lambda}\Ass\,\cF_\lambda=\Ass\,\cF.
$$
\item
Let $f:Y\to X$ be a flat morphism of schemes, and suppose that
the topological space $|X|$ underlying $X$ is locally noetherian.
Then $\Ass\,f^*\cF\subset f^{-1}\Ass\,\cF$.
\end{enumerate}
\end{proposition}
\begin{proof} (i) and (v) follow directly from the definitions.

(ii): Consider, for every point $x\in X$ the induced exact
sequence of $\cO_{\!X(x)}$-modules :
$$
0\to\underline\Gamma_{\{x\}}\cF'_x\to\underline\Gamma_{\{x\}}\cF_x
\to\underline\Gamma_{\{x\}}\cF''_x.
$$
Then :
$$
\Supp\,\underline\Gamma_{\{x\}}\cF'_x\subset
\Supp\,\underline\Gamma_{\{x\}}\cF_x\subset
\Supp\,\underline\Gamma_{\{x\}}\cF'_x\cup
\Supp\,\underline\Gamma_{\{x\}}\cF''_x
$$
which, in light of (i), is equivalent to the contention.

(iii): One applies (ii) to the natural injection :
$\cF/(\cG_1\cap\cdots\cap\cG_n)\to\oplus^n_{i=1}\cF/\cG_i$.

(iv): We may assume that $Y$ is a local scheme with
closed point $y\in Y$. Let $x_1,\dots,x_n\in X$ be the finitely
many points lying over $y$; for every $i,j\leq n$, we let
$$
\pi_j:X(x_j)\to Y\qquad \text{and} \qquad
\pi_{ij}:X(x_i)\times_X X(x_j)\to Y
$$
be the natural morphisms. To ease notation, denote also by
$\cF_i$ (resp. $\cF_{ij}$) the pull back of of $\cF$ to
$X(x_j)$ (resp. to $X(x_i)\times_X X(x_j)$). The induced morphism
$U:=X(x_1)\amalg\cdots\amalg X(x_n)\to X$ is faithfully flat,
so descent theory yields an exact sequence of $\cO_Y$-modules :
\set\begin{equation}\label{eq_descent}
\xymatrix{
0 \ar[r] & f_*\cF \ar[r] &
\prod^n_{j=1}\pi_{j*}\cF_j \ar@<.5ex>[r]^-{p_1^*}
\ar@<-.5ex>[r]_-{p_2^*} & \prod^n_{i,j=1}\pi_{ij*}\cF_{ij}.
}
\end{equation}
where $p_1,p_2:U\times_XU\to U$ are the natural morphisms.

\begin{claim}\label{cl_two-maps} The induced maps :
$$
\Gamma_{\!\{y\}}p_1^*,\Gamma_{\!\{y\}}p_2^*:
\prod^n_{j=1}\Gamma_{\!\{y\}}\pi_{j*}\cF_{_j}\to
\prod^n_{i,j=1}\Gamma_{\!\{y\}}\pi_{ij*}\cF_{ij}
$$
coincide.
\end{claim}
\begin{pfclaim} Indeed, it suffices to verify that they
coincide after projecting onto each factor
$\cG_{ij}:=\Gamma_{\!\{y\}}\pi_{ij*}\cF_{ij}$.
But this is clear from definitions if $i=j$. On the other
hand, if $i\neq j$, the image of $\pi_{ij}$ in $Y$ does
not contain $y$, so the corresponding factor $\cG_{ij}$
vanishes.
\end{pfclaim}

From \eqref{eq_descent} and claim \ref{cl_two-maps} we deduce
that the natural map
$$
\Gamma_{\!\{y\}}f_*\cF\to\prod^n_{j=1}\Gamma_{\!\{y\}}\pi_{j*}\cF_j
$$
is a bijection. Clearly
$\Gamma_{\!\{y\}}\pi_{j*}\cF_j=\Gamma_{\!\{x\}}\cF_j$, so both
assertions follow easily.

(vi): Let $x\in X$ be any point; since $|X|$ is locally noetherian,
the subset $\{x\}$ is closed and constructible in $X(x)$, and then
$f^{-1}(x)$ is closed and constructible in $Y\times_XX(x)$. Hence
$\delta(y,f^*\cF)\geq\delta(x,\cF)$ for every $y\in f^{-1}(x)$
(lemma \ref{lem_without-cohereur}(iii) and theorem
\ref{th_local-depth}). Then the assertion follows from (i).
\end{proof}

\begin{corollary}\label{cor_asso-tensor}
In the situation of corollary {\em\ref{cor_depth-flat-basechange}},
suppose furthermore that $|S|$ is a locally noetherian topological
space. Then we have :
$$
\Ass\,\cF\otimes_{\cO_{\!X}}f^*\cG=
\bigcup_{s\in\Ass\,\cG}\Ass\,\cF\otimes_{\cO_{S,s}}\kappa(s).
$$
\end{corollary}
\begin{proof} Under the stated assumptions, it is easily seen
that, for every $x\in X$ (resp. $s\in S$), the subset $\{x\}$
(resp. $\{s\}$) is constructible in $X(x)$ (resp. in $S(s)$).
Then the assertion follows immediately from proposition
\ref{prop_misc-Ass}(i) and theorem \ref{th_local-depth}.
\end{proof}

\begin{remark}
(i)\ \
Actually, it can be shown that the extra assumption on $|S|$
in corollary \ref{cor_asso-tensor} is superflous : see
\cite[Part I, Prop.3.4.3]{Gr-Ra}.

(ii)\ \
In the situation of proposition \ref{prop_misc-Ass}, one
could also consider the set of points $x\in X$ such that
$\delta'(x,\cF)=0$. This set contains $\Ass\,\cF$, by
proposition \ref{prop_misc-Ass}(i) and lemma
\ref{lem_ineq-deltas}(i). Such points are called
{\em attached primes} or {\em strong Krull primes} of $\cF$
in some literature : see {\em e.g.} \cite{Ep-Sha}, \cite{Dut}
(though the terminology ``attached prime'' is used also
for a different, unrelated notion: see the footnote on
page 404 of \cite{Dut}).
\end{remark}

\begin{corollary}\label{cor_prim-dec-local}
Let $X$ be a quasi-compact and quasi-separated scheme, $\cF$
a quasi-coherent $\cO_{\!X}$-module and $\cG\subset\cF$ a
quasi-coherent submodule. Then :

{\em (i)}\ \
For every quasi-compact open subset $U\subset X$ and every
quasi-coherent primary $\cO_{\!U}$-submodule $\cN\subset\cF_{|U}$,
with $\cG_{|U}\subset\cN$, there exists a quasi-coherent primary
$\cO_{\!X}$-submodule $\cM\subset\cF$ such that $\cM_{|U}=\cN$
and $\cG\subset\cN$.

{\em (ii)}\ \
$\cG$ admits a primary decomposition if and only if there
exists a finite open covering $X=U_1\cup\cdots\cup U_n$
consisting of quasi-compact open subsets, such that the
submodules $\cG_{|U_i}\subset\cF_{|U_i}$ admit primary
decompositions for every $i=1,\dots,n$.
\end{corollary}
\begin{proof} Say that $\cN$ is $x$-primary for some
point $x\in U$, and set $Z:=X\!\setminus\! U$. According to
\cite[Ch.I, Prop.9.4.2]{EGAI}, we can extend $\cN$ to a
quasi-coherent $\cO_{\!X}$-submodule $\cM_1\subset\cF$; up
to replacing $\cM_1$ by $\cM_1+\cG$, we may assume that
$\cG\subset\cM_1$. Since $(\cF/\cM_1)_{|U}=\cF_{|U}/\cN$,
it follows from lemma \ref{lem_Ass-Supp}(ii) that
$\Ass\,\cF/\cM_1\subset\{x\}\cup Z$. Let
$\bar\cM:=\underline\Gamma_Z(\cF/\cM_1)$, and denote by
$\cM$ the preimage of $\bar\cM$ in $\cF$. There follows
a short exact sequence :
$$
0\to\bar\cM\to\cF/\cM_1\to\cF/\cM\to 0.
$$
Clearly $R^1\underline\Gamma_Z\bar\cM=0$, whence a
short exact sequence
$$
0\to\underline\Gamma_Z\bar\cM\to\underline\Gamma_Z(\cF/\cM_1)
\to\underline\Gamma_Z(\cF/\cM)\to 0.
$$
We deduce that $\depth_Z(\cF/\cM)>0$; then theorem
\ref{th_local-depth} and proposition \ref{prop_misc-Ass}(i)
show that $Z\cap\Ass\,(\cF/\cM)=\emptyset$, so that $\cM$
is $x$-primary, as required.

(ii): We may assume that a covering $X=U_1\cup\cdots\cup U_n$
is given with the stated property. For every $i\leq n$,
let $\cG_{|U_i}=\cN_{i1}\cap\cdots\cap\cN_{ik_i}$ be a primary
decomposition; by (i) we may extend every $\cN_{ij}$ to
a primary submodule $\cM_{ij}\subset\cF$ containing
$\cG$. Then $\cG=\bigcap^n_{i=1}\bigcap^{k_i}_{j=1}\cM_{ij}$
is a primary decomposition of $\cG$.
\end{proof}

\sset\subsubsection{}\label{subsec_beaucoup-d-hypo}
Let $A$, $B$ be two local rings, with $A$ a domain,
$f:A\to B$ a local ring homomorphism, $M$ an $f$-flat
$B$-module of finite presentation. Denote by $K$ (resp.
$\kappa$) the field of fractions (resp. the residue
field) of $A$, and set
$$
B_0:=B\otimes_A\kappa
\qquad
B_K:=B\otimes_AK
\qquad
M_0:=M\otimes_A\kappa
\qquad
M_K:=M\otimes_AK.
$$

\begin{proposition} In the situation of
\eqref{subsec_beaucoup-d-hypo}, suppose that :
\begin{enumerate}
\alphaenu
\item
either, $A$ and $B$ are both noetherian rings
\item
or else, $f$ is essentially of finite presentation.
\end{enumerate}
Let $\fp\in\Ass_{B_K}M_K$ be any point, and denote by
$V(\fp)$ the topological closure of $\{\fp\}$
in $\Spec\,B$. Then all the maximal points of
$V(\fp)\cap\Spec\,B_0$ lie in $\Ass_{B_0}M_0$.
\end{proposition}
\begin{proof} We consider first case (a). Via the inclusion
map $\Spec\,B_K\to\Spec\,B$ we may regard $\fp$ as a prime
ideal of $B$, and $\Ass_{B_K}M_K$ as a subset of $\Ass_BM$;
by lemma \ref{lem_ass-are-mins}(ii), there exists a cyclic
$B$-submodule $N$ of $M$ such that $\Ann_BN=\fp$. Denote by
$\fm_A$ the maximal ideal of $A$, endow $M$ with its
$\fm_A$-adic filtration $\Fil^\bullet M$, and $N$ with the
filtration $\Fil^\bullet N$ induced by $\Fil^\bullet M$. Let
also $\gr^\bullet M$ and $\gr^\bullet N$ be the graded
$B_0$-modules associated with these filtrations on $M$ and $N$.
Since $M$ is $f$-flat, we have a $B_0$-linear isomorphism
$$
\gr^iM\isom(\fm^i_A/\fm^{i+1}_A)\otimes_\kappa M_0
\qquad
\text{for every $i\in\N$}
$$
which implies that $\Ass_B(\gr^iM)=\Ass_BM_0$ for every
$i\in\N$. From this, proposition \ref{prop_misc-Ass}(ii) and
an easy induction on $i\in\N$, we deduce that
$$
\Ass_B(N/\Fil^iN)\subset\Ass_B(M/\Fil^iM)\subset\Ass_BM_0
\qquad
\text{for every $i\in\N$}.
$$
On the other hand, by the Artin-Rees lemma (\cite[Th.8.5]{Mat}),
there exist $i,j\in\N$ such that
$$
\fm_A^jN\subset\Fil^iN\subset\fm_AN
$$
so that $\Supp\,N/\fm_AN=\Supp\,N/\Fil^iN$. However, lemma
\ref{lem_ass-are-mins}(i) implies that all the maximal
points of $\Supp\,N/\Fil^iN$ lie in $\Ass_B(N/\Fil^iN)$.
To conclude, it suffices to remark that the support of
$N/\fm_AN$ equals $V(\fp)\cap\Spec\,B_0$.

Next, suppose that (b) holds, and say that $B=C_\fq$ for
some finitely presented $A$-algebra $C$ and some prime
ideal $\fq\subset C$ whose preimage in $A$ is the maximal
ideal. We write $A$ as the colimit of a filtered system
$(A_\lambda~|~\lambda\in\Lambda)$ of noetherian local
subrings, such that the transition maps $A_\mu\to A_\lambda$
are local ring homomorphisms for every $\lambda,\mu\in\Lambda$
with $\lambda\geq\mu$; then there exist $\mu\in\Lambda$ and an
$A_\mu$-algebra $C_\mu$ of finite type with an isomorphism
of $C_\mu\otimes_{A_\mu}A\to C$ of $A$-algebras. For every
$\lambda\geq\mu$, let $C_\lambda:=C_\mu\otimes_{A_\mu}A_\lambda$,
denote by $\fq_\lambda\subset C_\lambda$ the preimage of $\fq$
under the induced map $C_\lambda\to C$, and set
$B_\lambda:=(C_\lambda)_{\fq_\lambda}$. Clearly $B$ is the
colimit of the resulting filtered system of rings
$(B_\lambda~|~\lambda\in\Lambda,\,\lambda\geq\mu)$, and
the induced maps $A_\lambda\to B_\lambda$ are local ring
homomorphisms, for every $\lambda\geq\mu$;
therefore we may assume that there exists a finitely
generated $B_\mu$-module $M_\mu$ with a $B$-linear
isomorphism $M_\mu\otimes_{B_\mu}B\isom M$, and we may
likewise define $M_\lambda:=M_\mu\otimes_{B_\mu}B_\lambda$
for every $\lambda\geq\mu$. In view of
\cite[Ch.IV, Cor.11.2.6.1(i)]{EGAIV-3}, we may further
assume that $M_\lambda$ is a flat $A_\lambda$-module, for
every $\lambda\geq\mu$, and after replacing $\Lambda$
by a cofinal subset, we may assume that $\mu$ is the
initial element of the partially ordered set $\Lambda$.
Then, to ease notation, denote by $K_\lambda$ (resp.
$\kappa_\lambda$) the field of fractions (resp. the
residue field) of $A_\lambda$, and set
$$
\begin{aligned}
B_{K,\lambda}:=\, & B_\lambda\otimes_{A_\lambda}K_\lambda
& \quad
M_{K,\lambda}:=\, & M_\lambda\otimes_{A_\lambda}K_\lambda
& \quad
C_{K,\lambda}:=\, & C_\lambda\otimes_{A_\lambda}K_\lambda \\
B_{0,\lambda}:=\, & B_\lambda\otimes_{A_\lambda}\kappa_\lambda
& \quad
M_{0,\lambda}:=\, & M_\lambda\otimes_{A_\lambda}\kappa_\lambda
\end{aligned}
\qquad
\text{for every $\lambda\in\Lambda$}.
$$
Notice that the resulting map
$f_\lambda:C_{K,\lambda}\to C_K:=C\otimes_AK$ induces
an isomorphism
$$
C_{K,\lambda}\otimes_{K_\lambda}K\isom C_K.
$$
So, $f_\lambda$ is a flat ring homomorphism, hence the
natural morphism $\phi_{K,\lambda}:\Spec\,B_K\to\Spec\,B_{K,\lambda}$
is also flat, and proposition \ref{prop_misc-Ass}(vi) yields
\set\begin{equation}\label{eq_asso-descend}
\Ass_{B_K}M_K\subset\phi^{-1}_{K,\lambda}\Ass_{B_{K,\lambda}}(M_{K,\lambda})
\qquad
\text{for every $\lambda\in\Lambda$}.
\end{equation}
Now, say that $\fp\in\Ass_{B_K}M_K$, and let
$\fp_\lambda\subset B_{K,\lambda}$ be the preimage of $\fp$,
for every $\lambda\in\Lambda$; we have
$\fp_\lambda\in\Ass_{B_{K,\lambda}}(M_{K,\lambda})$ by virtue
of \eqref{eq_asso-descend}, so in this situation the case
(a) of the proposition applies and shows that
$$
Z_\lambda:=\Max(\Spec\,B_{0,\lambda}/\fp_\lambda B_{0,\lambda})
\subset Z'_\lambda:=\Ass_{B_\lambda}(M_{0,\lambda})
\qquad
\text{for every $\lambda\in\Lambda$}
$$
(where, for a scheme $Z$, we have denoted $\Max(Z)$ the set
of maximal points of $Z$). Notice moreover, that the
natural morphism
$\phi_{0,\lambda}:\Spec\,B_0\to\Spec\,B_{0,\lambda}$
is also flat, so we have
$$
\Max(\Spec\,B_0/\fp_\lambda B_0)=
\bigcup_{z\in Z_\lambda}\Max(\phi^{-1}_{0,\lambda}(z))
\qquad
\text{for every $\lambda\in\Lambda$}
$$
by the going down theorem (\cite[Th.9.5]{Mat}). On the
other hand, from \cite[Ch.IV, Prop.4.2.7]{EGAIV-2} we have
$$
\bigcup_{z\in Z'_\lambda}\Max(\phi^{-1}_{0,\lambda}(z))=
(\Spec\,B_0)\cap\Ass_{B_{0,\lambda}\otimes_{k_\lambda}\kappa}
(M_{0,\lambda}\otimes_{\kappa_\lambda}\kappa)=\Ass_{B_0}(M_0)
$$
for every $\lambda\in\Lambda$. Summing up, we are reduced
to checking that there exists $\lambda\in\Lambda$ such that
the surjection $B_0/\fp_\lambda B_0\to B_0/\fp B_0$ is an
isomorphism, {\em i.e.} such that $\fp_\lambda B_0=\fp B_0$.
However, clearly $\fp B_0$ is the filtered union of the
system of subideals $(\fp_\lambda B_0~|~\lambda\in\Lambda)$;
but $B_0$ is a noetherian ring, so the claim is obvious.
\end{proof}

\begin{lemma}\label{lem_a&b}
Let $X$ be a quasi-separated and quasi-compact scheme, $\cF$
a quasi-coherent $\cO_{\!X}$-module of finite type. Then the
submodule $0\subset\cF$ admits a primary decomposition if
and only if the following conditions hold:
\begin{enumerate}
\alphaenu
\item
$\Ass\,\cF$ is a finite set.
\item
For every $x\in\Ass\,\cF$ there is a $x$-primary ideal
$I\subset\cO_{\!X,x}$ such that the natural map
$$
\Gamma_{\!\{x\}}\cF_x\to\cF_x/I\cF_x
$$
is injective.
\romanenu
\end{enumerate}
\end{lemma}
\begin{proof} In view of corollary \ref{cor_prim-dec-local}(i)
we are reduced to the case where $X$ is an affine scheme,
say $X=\Spec\,A$, and $\cF=M^\sim$ for an $A$-module $M$ of
finite type. Suppose first that $0$ admits a primary decomposition :
\set\begin{equation}\label{eq_prim-decomp}
0=\bigcap^k_{i=1}N_i.
\end{equation}
Then the natural map $M\to\oplus^k_{i=1}M/N_i$ is injective,
hence
\set\begin{equation}\label{eq_decomp-ass}
\Ass\,M\subset\Ass\,\bigoplus^k_{i=1}M/N_i\subset
\bigcup^k_{i=1}\Ass\,M/N_i
\end{equation}
by proposition \ref{prop_misc-Ass}(ii), and this shows that (a)
holds. Next, if $N_i$ and $N_j$ are $\fp$-primary
for the same prime ideal $\fp\subset A$, we may replace
both of them by their intersection (proposition \ref{prop_misc-Ass}(iii)).
Proceeding in this way, we achieve that the $N_i$ appearing
in \eqref{eq_prim-decomp} are primary submodules for pairwise
distinct prime ideals $\fp_1,\dots,\fp_k\subset A$. By
\eqref{eq_decomp-ass} we have $\Ass\,M\subset\{\fp_1,\dots,\fp_k\}$.
Suppose that $\fp_1\notin\Ass\,M$, and set
$Q:=\Ker\,(M\to\oplus^k_{i=2}M/N_i)$. For every $j>1$
we have $\Gamma_{\!\{\fp_j\}}(M/N_1)^\sim=0$, therefore :
$$
\Gamma_{\!\{\fp_j\}}Q^\sim_{\fp_j}=
\Ker\,(\Gamma_{\!\{\fp_j\}}M_{\fp_j}^\sim\to
\oplus^k_{i=1}(M/N_i)_{\fp_j}^\sim)=0.
$$
Hence $\Ass\,Q=\emptyset$, by proposition \ref{prop_misc-Ass}(i),
and then $Q=0$ by lemma \ref{lem_Ass-Supp}(iii). In other words,
we can omit $N_1$ from \eqref{eq_prim-decomp}, and still obtain
a primary decomposition of $0$; iterating this argument, we
may achieve that $\Ass\,M=\{\fp_1,\dots,\fp_k\}$. After these
reductions, we see that
$$
\Gamma_{\!\{\fp_j\}}(M/N_i)_{\fp_j}^\sim=0\qquad\text{whenever $i\neq j$}
$$
and consequently :
$$
\Ker\,(\phi_j:\Gamma_{\!\{\fp_j\}}M^\sim_{\fp_j}\to(M/N_j)_{\fp_j})=0
\qquad\text{for every $j=1,\dots,k$}.
$$
Now, for given $j\leq k$, let $\bar f_1,\dots,\bar f_n$
be a system of non-zero generators for $M/N_j$; by assumption
$I_i:=\Ann_A(\bar f_i)$ is a $\fp_j$-primary ideal for every
$i\leq n$. Hence $\fq_j:=\Ann_A(M/N_j)=\bigcap^n_{i=1}I_i$
is $\fp_j$-primary as well; since $\phi_j$ factors through
$(M/\fq_j M)_{\fp_j}$, we see that (b) holds.

Conversely, suppose that (a) and (b) hold. Say that
$\Ass\,M=\{\fp_1,\dots,\fp_k\}$; for every $j\leq k$
we choose a $\fp_j$-primary ideal $\fq_j$ such that
the map
$\Gamma_{\!\{\fp_j\}}M^\sim_{\fp_j}\to M_{\fp_j}/\fq_jM_{\fp_j}$
is injective. Clearly $N_j:=\Ker\,(M\to M_{\fp_j}/\fq_jM_{\fp_j})$
is a $\fp_j$-primary submodule of $M$; moreover, the induced map
$\phi:M\to\oplus^k_{j=1}M/N_j$ is injective, since
$\Ass\,\Ker\,\phi=\emptyset$ (proposition \ref{prop_misc-Ass}(ii)
and lemma \ref{lem_Ass-Supp}(iii)). In other words,
$0=\bigcap_{j=1}^kN_j$ is a primary decomposition of $0$.
\end{proof}

\begin{proposition}\label{prop_sh-exact_dec}
Let $X$ be a quasi-compact and quasi-separated scheme, and
$0\to\cF'\to\cF\to\cF''\to 0$ a short exact sequence of
quasi-coherent $\cO_{\!X}$-modules of finite type such that
\begin{enumerate}
\alphaenu
\item
$\Ass\,\cF'\cap\Supp\,\cF''=\emptyset$.
\item
The submodules $0\subset\cF'$ and $0\subset\cF''$ admit primary
decompositions.
\romanenu
\end{enumerate}
Then the submodule $0\subset\cF$ admits a primary decomposition.
\end{proposition}
\begin{proof} The assumptions imply that $\Ass\,\cF'$ and $\Ass\,\cF''$
are finite sets, (lemma \ref{lem_a&b}), hence the same holds for
$\Ass\,\cF$ (proposition \ref{prop_misc-Ass}(ii)). Given any point
$x\in X$ and any $x$-primary ideal $I\subset\cO_{\!X,x}$, we may
consider the commutative diagram with exact rows:
$$
\xymatrix{
0 \ar[r] & \Gamma_{\!\{x\}}\cF'_x \ar[r] \ar[d]_\alpha &
\Gamma_{\!\{x\}}\cF_x \ar[r] \ar[d]_\beta &
\Gamma_{\!\{x\}}\cF''_x \ar[d]_\gamma \\
\Tor_1^{\cO_{\!X,x}}(I,\cF''_x) \ar[r] & \cF'_x/I\cF'_x \ar[r] &
\cF_x/I\cF_x \ar[r] & \cF''_x/I\cF''_x
}
$$
Now, if $x\in\Ass\,\cF'$, assumption (a) implies that
$\Gamma_{\!\{x\}}\cF''_x=0=\Tor_1^{\cO_{\!X,x}}(I_x,\cF''_x)$, and by
(b) and lemma \ref{lem_a&b} we can choose $I$ such that $\alpha$ is injective;
a little diagram chase then shows that $\beta$ is injective as well.
Similarly, if $x\in\Ass\,\cF''$, we have $\Gamma_{\!\{x\}}\cF'_x=0$
and we may choose $I$ such that $\gamma$ is injective, which implies
again that $\beta$ is injective. To conclude the proof it suffices to
apply again lemma \ref{lem_a&b}.
\end{proof}

\begin{proposition}\label{prop_prim-dec&fin-map}
Let $Y$ be a quasi-compact and quasi-separated scheme, $f:X\to Y$
a finite morphism and $\cF$ a quasi-coherent $\cO_{\!X}$-module of
finite type. Then the $\cO_{\!X}$-submodule $0\subset\cF$ admits a
primary decomposition if and only if the same holds for the
$\cO_Y$-submodule $0\subset f_*\cF$.
\end{proposition}
\begin{proof} Under the stated assumptions, we can apply
the criterion of lemma \ref{lem_a&b}. To start out, it is
clear from proposition \ref{prop_misc-Ass}(iv.b) that $\Ass\,\cF$
is finite if and only if $\Ass\,f_*\cF$ is finite.
Next, suppose that $0\subset\cF$ admits a primary decomposition,
let $y\in\Ass\,f_*\cF$ be any point and set
$f^{-1}(y):=\{x_1,\dots,x_n\}$; for every $j\leq n$ we can find
an $x_j$-primary ideal $I_j\subset\cO_{\!X}$ such that the map
$\Gamma_{\!\{x_j\}}\cF_{x_j}\to\cF_{x_j}/I_j\cF_{x_j}$
is injective. Let $I$ be the kernel of the natural map
$\cO_Y\to f_*(\cO_{\!X}/(I_1\cap\cdots\cap I_n))$.
Then $I$ is a quasi-coherent $y$-primary ideal of $\cO_Y$
and we deduce a commutative diagram :
\set\begin{equation}\label{eq_diag-up&down}
{\diagram
\bigoplus^n_{j=1}\Gamma_{\!\{x_j\}}(\cF_{x_j}) \ar[r] \ar[d]_\alpha &
\Gamma_{\!\{y\}}(f_*\cF_y) \ar[d]^\beta \\
\bigoplus^n_{j=1}\cF_{x_j}/I\cF_{x_j} \ar[r] & f_*\cF_y/If_*\cF_y
\enddiagram}\end{equation}
whose horizontal arrows are isomorphisms, in view of
proposition \ref{prop_misc-Ass}(iv.a), and where $\alpha$ is
injective by construction. It follows that $\beta$ is
injective, so condition (b) of lemma \ref{lem_a&b} holds
for the stalk $f_*\cF_y$, and since $y\in\Ass\,f_*\cF$ is
arbitrary, we see that $0\subset f_*\cF$ admits a primary
decomposition. Conversely, suppose that $0\subset f_*\cF$
admits a primary decomposition; then for every $y\in Y$ we
can find a quasi-coherent $y$-primary ideal $I\subset\cO_Y$
such that the map $\beta$ of \eqref{eq_diag-up&down} is
injective; hence $\alpha$ is injective as well, and again
we deduce easily that $0\subset\cF$ admits a primary
decomposition.
\end{proof}

\begin{proposition}\label{prop_prim_op&clo}
Let $X$ and $\cF$ be as in lemma {\em\ref{lem_a&b}}, $i:Z\to X$
a closed constructible immersion, and $U:=X\!\setminus\! Z$.
Suppose that :
\begin{enumerate}
\alphaenu
\item
The $\cO_{\!U}$-submodule $0\subset\cF_{|U}$ and the $\cO_{\!Z}$-submodule
$0\subset i^*\cF$ admit primary decompositions.
\item
The natural map\ \ $\underline\Gamma_Z\cF\to i_*i^*\cF$\ \ is injective.
\romanenu
\end{enumerate}
Then the $\cO_{\!X}$-submodule $0\subset\cF$ admits a primary decomposition.
\end{proposition}
\begin{proof} We shall verify that conditions (a) and (b) of lemma
\ref{lem_a&b} hold for $\cF$. To check condition (a), it suffices
to remark that $U\cap\Ass\,\cF=\Ass\,\cF_{|U}$ (which is obvious)
and that $Z\cap\Ass\,\cF\subset\Ass\,i^*\cF$, which follows easily
from our current assumption (b).

Next we check that condition (b) of {\em loc.cit.} holds. This
is no problem for the points $x\in U\cap\Ass\,\cF$, so suppose that
$x\in Z\cap\Ass\,\cF$. Moreover, we may also assume that $X$ is
affine, say $X=\Spec\,A$, so that $Z=V(J)$ for some ideal $J\subset A$.
Due to proposition \ref{prop_prim-dec&fin-map} we know that
$0\subset\cG:=i_*i^*\cF$ admits a primary decomposition, hence we can
find an $x$-primary ideal $I\subset A$ such that the natural map
$\Gamma_{\!\{x\}}\cG\to\cG_x/I\cG_x\simeq\cF_x/(I+J)\cF_x$ is injective.
Clearly $I+J$ is again a $x$-primary ideal; since $Z$ is closed
and constructible, corollary \ref{cor_better-late} and our assumption
(b) imply that the natural map $\Gamma_{\!\{x\}}\cF\to\Gamma_{\!\{x\}}\cG$
is injective, hence the same holds for the map
$\Gamma_{\!\{x\}}\cF\to\cF_x/(I+J)\cF_x$, as required.
\end{proof}

\sset\subsubsection{}
Let now $A$ be a ring, $\fp\subset A$ any prime ideal, and
$n\geq 1$ any integer; if $A/\fp^n$ does not admit imbedded
primes, then the ideal $\fp^n$ is $\fp$-primary. In the
presence of imbedded primes, this fails. For instance, we
have the following :

\begin{example} Let $k$ be a field, $k[x,y]$ the free polynomial
algebra in indeterminates $x$ and $y$; consider the ideal
$I:=(xy,y^2)$, and set $A:=k[x,y]/I$. Let $\bar x$ and $\bar y$
be the images of $x$ and $y$ in $A$; then $\fp:=(\bar y)\subset A$
is a prime ideal, and $\fp^2=0$. However, $\fm:=\Ann_A(\bar y)$
is the maximal ideal generated by $\bar x$ and $\bar y$, so
the ideal $0\subset A$ is not $\fp$-primary.
\end{example}

There is however a natural sequence of $\fp$-primary ideals
naturally attached to $\fp$. To explain this, let us remark,
more generally, the following :

\begin{lemma}\label{lem_bije-primary}
Let $A$ be a ring, $\fp\subset A$ a prime ideal. Denote by
$\phi_\fp:A\to A_\fp$\ \ the localization map. The rule :
$$
I\mapsto\phi^{-1}_\fp I
$$
induces a bijection from the set of $\fp A_\fp$-primary
ideals of $A_\fp$, to the set of $\fp$-primary ideals of $A$.
\end{lemma}
\begin{proof} Suppose that $I\subset A_\fp$ is $\fp A_\fp$-primary;
since the natural map $A/\phi_\fp^{-1}I\to A_\fp/I$ is injective,
it is clear that $\phi_\fp^{-1}I$ is $\fp$-primary. Conversely,
suppose that $J\subset A$ is $\fp$-primary; we claim that
$J=\phi_\fp^{-1}(J_\fp)$. Indeed, by assumption (and by lemma
\ref{lem_Ass-Supp}(iii)) we have $(A/J)_\fq=0$ whenever
$\fq\neq\fp$; it follows easily that the localization map
$A/J\to A_\fp/J_\fp$ is an isomorphism, whence the contention.
\end{proof}

\begin{definition}
Keep the notation of lemma \ref{lem_bije-primary}; for every
$n\geq 0$ one defines the {\em $n$-th symbolic power\/} of
$\fp$, as the ideal :
$$
\fp^{(n)}:=\phi^{-1}_\fp(\fp^nA_\fp).
$$
By lemma \ref{lem_bije-primary}, the ideal $\fp^{(n)}$ is
$\fp$-primary for every $n\geq 1$. More generally, for every
$A$-module $M$, let $\phi_{M,\fp}:M\to M_\fp$ be the localization
map; then one defines the {\em $\fp$-symbolic filtration\/} on $M$,
by the rule :
$$
\Fil_\fp^{(n)}M:=\phi_{M,\fp}^{-1}(\fp^nM_\fp)
\qquad\text{for every $n\geq 0$}.
$$
The filtration $\Fil^{(\bullet)}_\fp M$ induces a (linear)
topology on $M$, called the {\em $\fp$-symbolic topology}.

More generally, if $\Sigma\subset\Spec\,A$ is any subset,
we define the {\em $\Sigma$-symbolic topology\/} on $M$,
as the coarsest linear topology $\cT_\Sigma$ on $M$ such
that $\Fil_\fp^{(n)}M$ is an open subset of $\cT_\Sigma$,
for every $\fp\in\Sigma$ and every $n\geq 0$. If $\Sigma$
is finite, it is induced by the {\em $\Sigma$-symbolic filtration},
defined by the submodules :
$$
\Fil^{(n)}_\Sigma M:=\bigcap_{\fp\in\Sigma}\Fil^{(n)}_\fp M
\qquad\text{for every $n\geq 0$}.
$$
\end{definition}

\sset\subsubsection{}\label{subsec_dubbydubby}
Let $A$ be a ring, $M$ an $A$-module, $I\subset A$ an ideal.
We shall show how to characterize the finite subsets
$\Sigma\subset\Spec\,A$ such that the $\Sigma$-symbolic
topology on $M$ agrees with the $I$-preadic topology (see
theorem \ref{th_brain-dead}). Hereafter, we begin with
a few preliminary observations. First, suppose that $A$
is noetherian; then it is easily seen that, for every
prime ideal $\fp\subset A$, and every $\fp A_\fp$-primary
ideal $I\subset A_\fp$, there exists $n\in\N$ such that
$\fp^nA_\fp\subset I$. From lemma \ref{lem_bije-primary}
we deduce that every $\fp$-primary ideal of $A$ contains
a symbolic power of $\fp$; {\em i.e.}, every $\fp$-primary ideal
is open in the $\fp$-symbolic topology of $A$. More generally,
let $\Sigma\subset\Spec\,A$ be any subset, $M$ a finitely
generated $A$-module, and $N\subset M$ a submodule; from
the existence of a primary decomposition for $N$
(\cite[Th.6.8]{Mat}), we see that $N$ is open in the
$\Sigma$-symbolic topology of $M$, whenever
$\Ass\,M/N\subset\Sigma$. Especially, if
$\Ass\,M/I^nM\subset\Sigma$ for every $n\in\N$, then the
$\Sigma$-symbolic topology is finer than the $I$-preadic
topology on $M$.
On the other hand, clearly $I^nM\subset\Fil^{(n)}_\fp M$ for
every $\fp\in\Supp\,M/IM$ and every $n\in\N$, so the
$I$-preadic topology on $M$ is finer than the $\Supp\,M/IM$-symbolic
topology. Summing up, if we have :
$$
\Sigma_0(M):=\bigcup_{n\in\N}\Ass\,M/I^nM\subset\Sigma\subset\Supp\,M/IM
$$
then the $\Sigma$-symbolic topology on $M$ agrees with the
$I$-preadic topology. Notice that the above expression for
$\Sigma_0(M)$ is a union of finite subsets (\cite[Th.6.5(i)]{Mat}),
hence $\Sigma_0(M)$ is -- {\em a priori} -- at most countable; in
fact, we shall show that $\Sigma_0(M)$ is finite. Indeed, for every
$n\in\N$, set :
$$
\gr_I^nA:=I^n/I^{n+1}
\qquad
\gr^n_IM:=I^nM/I^{n+1}M.
$$
Then $\gr^\bullet_IA:=\oplus_{n\in\N}\gr^n_IA$ is naturally
a graded $A/I$-algebra, and $\gr^\bullet_IM:=\oplus_{n\in\N}\gr^n_IM$
is a graded $\gr^\bullet_IA$-module. Let
$\psi:\Spec\,\gr^\bullet_IA\to\Spec\,A/I$ be the natural
morphism, and set :
$$
\Sigma(M):=\psi(\Ass_{\gr^\bullet_IA}(\gr^\bullet_IM)).
$$

\begin{lemma}\label{lem_bobbybobby}
With the notation of \eqref{subsec_dubbydubby}, we have :
$$
\Ass_{A/I}(\gr_I^nM)\subset\Sigma(M)
\qquad\text{for every $n\in\N$}.
$$
\end{lemma}
\begin{proof}To ease notation, set $A_0:=A/I$ and
$B:=\gr_I^\bullet A$. Suppose that $\fp\in\Ass_{A_0}(\gr_I^nM)$;
by lemma \ref{lem_ass-are-mins}(i), there exists $m\in\gr_I^nM$
such that $\fp$ is the preimage of a minimal prime ideal of
$A_0/\Ann_{A_0}(m)$.
However, if we regard $m$ as a homogeneous element of the
$B$-module $\gr_I^\bullet M$, we have the obvious identity :
$$
\Ann_{A_0}(m)=A_0\cap\Ann_B(m)\subset B
$$
from which we see that the induced map
$$
A_0/\Ann_{A_0}(m)\to B/\Ann_B(m)
$$
is injective, hence $\psi$ restricts to a dominant morphism
$V(\Ann_B(m))\to V(\Ann_{A_0}(m))$. Especially, there exists
$\fq\in V(\Ann_B(m))$ with $\psi(\fq)=\fp$; up to replacing
$\fq$ by a generization, we may assume that $\fq$ is a maximal
point of $V(\Ann_B(m))$, hence $\fq$ is an associated prime
for $\gr_I^\bullet M$, again by lemma \ref{lem_ass-are-mins}(i).
\end{proof}

\sset\subsubsection{}\label{subsec_sigma-min}
An easy induction, starting from lemma \ref{lem_bobbybobby},
shows that $\Ass_AM/I^nM\subset\Sigma(M)$, for every $n\in\N$,
therefore $\Sigma_0(M)\subset\Sigma(M)$. However, if $A$ is noetherian,
the same holds for $\gr^\bullet_IA$ (since the latter is a quotient
of an $A/I$-algebra of finite type), hence $\Sigma(M)$ is finite,
provided $M$ is finitely generated (\cite[Th.6.5(i)]{Mat}), and
{\em a fortiori}, the same holds for $\Sigma_0(M)$.

\begin{remark} Another proof of the finiteness of $\Sigma_0(M)$
can be found in \cite{Brod}.
\end{remark}

\sset\subsubsection{}
Next, we wish to show that actually there exists a {\em smallest\/}
subset $\Sigma_{\min}(M)\subset\Spec\,A$ such that the
$\Sigma_{\min}(M)$-symbolic topology on $M$ agrees with $I$-preadic
topology; after some simple reductions, this boils down to the
following assertion.
Let $\Sigma\subset\Spec\,A$ be a finite subset, and $\fp,\fp'\in\Sigma$
two elements, such that the $\Sigma$-symbolic topology on $M$ agrees
with both the $\Sigma\setminus\!\{\fp\}$-symbolic topology and
the $\Sigma\!\setminus\!\{\fp'\}$-symbolic topology; then these
topologies agree as well with the
$\Sigma\setminus\!\{\fp,\fp'\}$-symbolic topology.
Indeed, given any subset $\Sigma'\subset\Spec\,A$ with $\fp\in\Sigma'$,
for the $\Sigma'$-symbolic and the
$\Sigma'\setminus\!\{\fp\}$-symbolic topologies to agree, it is
necessary and sufficient that, for every $n\in\N$ there exists
$m\in\N$ such that :
$$
\Fil_{\Sigma'\setminus\!\{\fp\}}^{(m)}M
\subset\Fil_\fp^{(n)}M
\qquad\text{or, what is the same :}\qquad
\Fil_{\Sigma'\setminus\!\{\fp\}}^{(m)}M_\fp
\subset\fp^nM_\fp.
$$
In the latter inclusion, we may then replace $\Sigma'\setminus\!\{\fp\}$
by the smaller subset $\Sigma'\cap\Spec\,A_\fp\setminus\!\{\fp\}$,
without changing the two terms. Suppose now that $\fp'\notin\Spec\,A_\fp$;
then if we apply the above, first with $\Sigma':=\Sigma$, and then
with $\Sigma':=\Sigma\setminus\!\{\fp'\}$, we see that the $\Sigma$-symbolic
topology agrees with the $\Sigma\setminus\!\{\fp\}$-symbolic topology,
if and only if the $\Sigma\setminus\!\{\fp'\}$-symbolic topology
agrees with the $\Sigma\setminus\!\{\fp,\fp'\}$-symbolic topology,
whence the contention. In case $\fp'\in\Spec\,A_\fp$, we may assume
that $\fp\neq\fp'$, otherwise there is nothing to prove; then we
shall have $\fp\notin\Spec\,A_{\fp'}$, so the foregoing argument
still goes through, after reversing the roles of $\fp$ and $\fp'$.

\sset\subsubsection{}
Finally, theorem \ref{th_brain-dead} will characterize
the subset $\Sigma_{\min}(M)$ as in \eqref{subsec_sigma-min}.
To this aim, for every prime ideal $\fp\subset A$, let $A^\wedge_\fp$
denote the $\fp$-adic completion of $A$; we set :
$$
\Ass_A(I,M):=\{\fp\in\Spec\,A~|~
\text{there exists
$\fq\in\Ass_{A^\wedge_\fp}M\!\otimes_A\!A^\wedge_\fp$\ \
such that\ \
$\textstyle\sqrt{\fq+IA^\wedge_\fp}=\fp A^\wedge_\fp$}\}
$$
(where, for any ideal $J\subset A^\wedge_\fp$ we denote by
$\sqrt{J}\subset A^\wedge_\fp$ the radical of $J$, so the above
condition selects the points $\fq\in\Spec\,A^\wedge_\fp$, such
that $\overline{\{\fq\}}\cap V(IA^\wedge_\fp)=\{\fp A^\wedge_\fp\}$).

\begin{lemma}\label{lem_little-used}
Let $A$ be a noetherian ring, $M$ an $A$-module. Then we have :
\begin{enumerate}
\item
$\depth_{A_\fp}M\!\otimes_A\!A_\fp=
\depth_{A^\wedge_\fp}M\!\otimes_A\!A^\wedge_\fp$\ \  for every
$\fp\in\Spec\,A$.
\item
$\Ass_A(0,M)=\Ass_AM$.
\item
$\Ass_A(I,M)\subset V(I)\cap\Supp\,M$.
\item
Suppose that $M$ is the union of a filtered family
$(M_\lambda~|~\lambda\in\Lambda)$ of $A$-submodules.
Then :
$$
\Ass_A(I,M)=\bigcup_{\lambda\in\Lambda}\Ass_A(I,M_\lambda).
$$
\item
$\Ass_A(I,M)$ contains the maximal points of\/ $V(I)\cap\Supp\,M$.
\end{enumerate}
\end{lemma}
\begin{proof} (iii) is immediate, and in view of \cite[Th.8.8]{Mat},
(i) is a special case of corollary \ref{cor_f-flat-invariance}.

(ii): By definition, $\Ass_A(0,M)$ consists of all the prime
ideals $\fp\subset A$ such that $\fp A_\fp^\wedge$ is an
associated prime ideal of $M\!\otimes_A\!A_\fp^\wedge$; in light
of proposition \ref{prop_misc-Ass}(i), the latter condition
holds if and only if
$\depth_{A^\wedge_\fp}M\!\otimes_A\!A^\wedge_\fp=0$, hence
if and only if $\depth_{A_\fp}M\!\otimes_A\!A_\fp=0$, by (i);
to conclude, one appeals again to proposition \ref{prop_misc-Ass}(i).

(iv): In view of proposition \ref{prop_misc-Ass}(v), we have
$\Ass_{A_\fp^\wedge}M\!\otimes_A\!A^\wedge_\fp=\bigcup_{\lambda\in\Lambda}
\Ass_{A_\fp^\wedge}M_\lambda\!\otimes_A\!A^\wedge_\fp$; the contention
is an immediate consequence.

(v): In view of (iv), we are easily reduced to the case where
$M$ is a finitely generated $A$-module, in which case $\Supp\,M=V(J)$
for some ideal $J\subset A$. Hence, suppose that $\fp\subset A$
is a maximal point of $\Spec\,A/(I+J)$, in other words, the preimage
of a minimal prime ideal of $A/(I+J)$; notice that
$\Supp\,M\!\otimes_A\!A^\wedge_\fp=V(JA^\wedge_\fp)$.
Hence we have $(J+I)A^\wedge_\fp\subset\fq+IA^\wedge_\fp$ for
every $\fq\in\Supp\,M\!\otimes_A\!A^\wedge_\fp$, so the radical
of $\fq+IA^\wedge_\fp$ equals $\fp A^\wedge_\fp$, as required.
\end{proof}

\begin{lemma}\label{lem_chevalley}
Let $A$ be a local noetherian ring, $\fm$ its maximal ideal,
and suppose that $A$ is $\fm$-adically complete.
Let $M$ be a finitely generated $A$-module, and $(\Fil^nM~|~n\in\N)$
a descending filtration consisting of $A$-submodules of $M$.
Then the following conditions are equivalent :
\begin{enumerate}
\alphaenu
\item
$\bigcap_{n\in\N}\Fil^nM=0$.
\item
For every $i\in\N$ there exists $n\in\N$ such that
$\Fil^nM\subset\fm^iM$.
\end{enumerate}
\end{lemma}
\begin{proof} Clearly (b)$\Rightarrow$(a), hence suppose that
(a) holds. For every $i,n\in\N$, set
$$
J_{i,n}:=\Img(\Fil^nM\to M/\fm^iM).
$$
For given $i\in\N$, the $A$-module $M/\fm^iM$ is artinian,
hence there exists $n\in\N$ such that $J_i:=J_{i,n}=J_{i,n'}$
for every $n'\geq n$. Assertion (b) then follows from the following :

\begin{claim} If (a) holds, $J_i=0$ for every $i\in\N$.
\end{claim}
\begin{pfclaim}[] By inspecting the definition, it is easily seen
that the natural $A/\fm^{i+1}$-linear map $J_{i+1}\to J_i$ is
surjective for every $i\in\N$, hence we are reduced to showing
that $J:=\lim_{i\in\N}\,J_i$ vanishes. However, $J$ is naturally
a submodule of $\lim_{i\in\N}\,M/\fm^iM\simeq M$, and if
$x\in M$ lies in $J$, then we have $x\in\Fil^n M+\fm^iM$ for every
$i,n\in\N$. Since $\Fil^nM$ is a closed subset for the $\fm$-adic
topology of $M$ (\cite[Th.8.10(i)]{Mat}), we have
$\bigcap_{i\in\N}(\Fil^n M+\fm^iM)=\Fil^nM$, for every $n\in\N$,
hence $x\in\bigcap_{n\in\N}\Fil^nM$, which vanishes, if (a) holds.
\end{pfclaim}
\end{proof}

The following theorem generalizes \cite[Ch.VIII, \S5, Cor.5]{Z-S}
and \cite[Prop.7.1]{Har2}.

\begin{theorem}\label{th_brain-dead}
Let $A$ be a noetherian ring, $I\subset A$ an ideal, $M$
a finitely generated $A$-module, and $\Sigma\subset\Spec\,A/I$
a finite subset. Then the $\Sigma$-symbolic topology on $M$
agrees with the $I$-preadic topology if and only if
$\Ass_A(I,M)\subset\Sigma$.
\end{theorem}
\begin{proof} Let $\fp\in\Ass_A(I,M)\!\setminus\!\Sigma$, and
suppose, by way of contradiction, that the $\Sigma$-symbolic
topology agrees with the $I$-adic topology, {\em i.e.} for every
$n\in\N$ there exists $m\in\N$ such that :
$$
\Fil^{(m)}_\Sigma M\subset I^nM.
$$
Let $X:=\Spec\,A_\fp$ and $U:=X\!\setminus\!\{\fp\}$; after
localizing at the prime $\fp$, we deduce the inclusion :
\set\begin{equation}\label{eq_loc-at-fp}
\Fil^{(m)}_{\Sigma\cap U} M_\fp\subset I^nM_\fp
\end{equation}
(cp. the discussion in \eqref{subsec_sigma-min}). Let also
$\cM:=M^\sim$, the quasi-coherent $\cO_{\!X}$-module associated
to $M$; clearly we have $I^mM\subset\Fil^{(m)}_\fq M$ for every
$\fq\in\Spec\,A/I$, hence \eqref{eq_loc-at-fp} implies the inclusion :
$$
\{x\in M_\fp~|~x_{|U}\in I^m\!\cM(U)\}\subset I^nM_\fp.
$$
Let $A^\wedge_\fp$ (resp. $M^\wedge_\fp$) be the $\fp$-adic
completion of $A$ (resp. of $M$), $f:X^\wedge:=\Spec\,A^\wedge_\fp\to X$
the natural morphism, and set
$U^\wedge:=f^{-1}U=X^\wedge\setminus\{\fp A^\wedge_\fp\}$.
Since $f$ is faithfully flat, and $U$ is quasi-compact, we deduce that :
\set\begin{equation}\label{eq_deduced}
\{x\in M^\wedge_\fp~|~x_{|U^\wedge}\in I^mf^*\!\cM(U^\wedge)\}
\subset I^nM^\wedge_\fp.
\end{equation}
On the other hand, by assumption there exists
$\fq\in\Ass_{A^\wedge_\fp}M^\wedge_\fp$ such that
$\overline{\{\fq\}}\cap V(IA^\wedge_\fp)=\{\fp A^\wedge_\fp\}$,
hence we may find $x\in M^\wedge_\fp$ whose support is
$\overline{\{\fq\}}$ (lemma \ref{lem_ass-are-mins}(ii)).
It follows easily that the image of $x$ vanishes in
$f^*\!\cM/I^mf^*\!\cM(U)$ for all $m\in\N$, {\em i.e.}
$x_{|U}\in I^mf^*\!\cM(U)$ for every $m\in\N$, hence
$x\in I^nM_\fp^\wedge$ for every $n\in\N$, in view of
\eqref{eq_deduced}. However, $M^\wedge_\fp$ is separated
for the $\fp$-adic topology, {\em a fortiori\/} also for
the $I$-preadic topology, so $x=0$, a contradiction.

Next, let $\fp\in\Sigma\!\setminus\!\Ass_A(I,M)$, set
$\Sigma':=\Sigma\!\setminus\!\{\fp\}$, and suppose that
the $\Sigma$-symbolic topology agrees with the $I$-adic
topology; we have to prove that the latter agrees as well
with the $\Sigma'$-symbolic topology. This amounts to
showing that, for every $n\in\N$ there exists $m\in\N$ such that :
$$
\Fil^{(m)}_{\Sigma'}M\subset\Fil^{(n)}_\fp M
\qquad\text{or, what is the same :}\qquad
\Fil^{(m)}_{\Sigma'}M_\fp\subset\fp^n M_\fp.
$$
We may write :
$$
\Fil^{(m)}_{\Sigma'}M_\fp=
\{x\in M_\fp~|~x_{|U}\in(\Fil^{(m)}_{\Sigma'}M_\fp)^\sim(U)\}=
\{x\in M_\fp~|~x_{|U}\in(\Fil^{(m)}_ \Sigma M_\fp)^\sim(U)\}
$$
from which it is clear that the contention holds if and only
if, for every $n\in\N$ there exists $m\in\N$ such that :
$$
\{x\in M_\fp~|~x_{|U}\in I^m\!\cM(U)\}\subset\fp^nM_\fp.
$$
Arguing as in the foregoing, we see that the latter condition
holds if and only if :
$$
\Fil^mM^\wedge_\fp:=
\{x\in M^\wedge_\fp~|~x_{|U^\wedge}\in I^mf^*\!\cM(U^\wedge)\}
\subset\fp^nM^\wedge_\fp.
$$
In view of lemma \ref{lem_chevalley}, we are then reduced to
showing the following :

\begin{claim} $\bigcap_{m\in\N}\Fil^mM^\wedge_\fp=0$.
\end{claim}
\begin{pfclaim}[] Let $x$ be an element in this intersection;
for any $\fq\in U^\wedge\cap V(IA^\wedge_\fp)$, we have
$x_\fq\in\bigcap_{m\in\N}\fq^m(M^\wedge_\fp)_\fq$, hence
$x_\fq=0$, by \cite[Th.8.10(i)]{Mat}. In other words,
$\Supp(x)\cap U\cap V(I)=\emptyset$. Suppose that $x\neq 0$,
and let $\fq$ be any maximal point of $\Supp(x)$;
by lemma \ref{lem_ass-are-mins}(i), $\fq$ is an associated
prime, and the foregoing implies that
$\{\fq\}\cap V(IA^\wedge_\fp)=\{\fp A^\wedge_\fp\}$, which
contradicts the assumption that $\fp\notin\Ass_A(I,M)$.
\end{pfclaim}
\end{proof}

\begin{corollary} In the situation of theorem
{\em\ref{th_brain-dead}}, $\Ass_A(I,M)$ is a finite set.
\end{corollary}
\begin{proof} We have already found a finite subset
$\Sigma\subset\Spec\,A/I$ such that the $\Sigma$-symbolic
topology on $M$ agrees with the $I$-preadic topology
(see \eqref{subsec_sigma-min}).
The contention then follows straightforwardly from
theorem \ref{th_brain-dead}.
\end{proof}

\begin{example} Let $k$ be a field with $\chara\,k\neq 2$,
and let $C\subset\A^2_k:=\Spec\,k[X,Y]$ be the nodal curve
cut by the equation $Y^2=X^2+X^3$, so that the only singularity
of $C$ is the node at the origin $p:=(0,0)\in C$. Let $R:=k[X,Y,Z]$,
$A:=R/(Y^2-X^2-X^3)$; denote by $\pi:\A^3_k:=\Spec\,R\to\A^2_k$ the
linear projection which is dual to the inclusion $k[X,Y]\to R$,
so that $D:=\pi^{-1}C=\Spec\,A$.
We define a morphism $\phi:\A^1_k:=\Spec\,k[T]\to D$ by the rule :
$T\mapsto(T^2-1,T(T^2-1),T)$ ({\em i.e.} $\phi$ is dual to the
homomorphism of $k$-algebras such that $X\mapsto T^2-1$,
$Y\mapsto T(T^2-1)$ and $Z\mapsto T$). Let $C'\subset D$ be
the image of $\phi$, with its reduced subscheme structure.
It is easy to check that the restriction of $\pi$ maps $C'$
birationally onto $C$, so there are precisely two points
$p'_0,p'_1\in C'$ lying over $p$. Let $\fn:=I(C')\subset A$,
the prime ideal which is the generic point of the (irreducible)
curve $C'$. We claim that the $\fn$-preadic topology on $A$
does not agree with the $\fn$-symbolic topology. To this aim
-- in view of theorem \ref{th_brain-dead} -- it suffices to
show that $\{p'_0,p'_1\}\subset\Ass_A(\fn,A)$. However, for any closed
point $\fp\in\pi^{-1}(p)$, the $\fp$-adic completion $A^\wedge_\fp$
admits two distinct minimal primes, corresponding to the two
branches of the nodal conic $C$ at the node $p$, and the
corresponding irreducible components of $B:=\Spec\,A^\wedge_\fp$
meet along the affine line $V(Z)$. To see this, we may suppose
that $\fp=(X,Y,Z)$, hence
$A^\wedge_\fp\simeq k[[X,Y,Z]]/(Y^2-X^2(1+X))$, and notice that
the latter is isomorphic to $k[[S,Y,Z]]/(Y^2-S^2)$, via the
isomorphism that sends $Y\mapsto Y$, $Z\mapsto Z$ and
$S\mapsto X(1+X)^{1/2}$ (the assumption on the characteristic
of $k$ ensures that $1+X$ admits a square root in $k[[X]]$).
Now, say that $\fp=p'_0$; then $C'_\fp:=C'\cap B$ is contained
in only one of the two irreducible components of $B$.
Let $\fq\in B$ be the minimal prime ideal whose closure does
not contain $C'_\fp$; then $\fq\in\Ass A^\wedge_\fp$ and
$\overline{\{\fq\}}\cap C'_\fp=\{\fp A^\wedge_\fp\}$, therefore
$p'_0\in\Ass_A(\fn,A)$, as stated.
\end{example}

\begin{example}\label{ex_Ass}
Let $A$ be an excellent normal ring, $I\subset A$ any ideal, and
set $Z:=V(I)$. Then $\Ass_A(I,A)$ is the set $\Max(Z)$ of
all maximal points of $Z$. Indeed, $\Max(Z)\subset\Ass_A(I,A)$
by lemma \ref{lem_little-used}(v).
Conversely, suppose $\fp\in\Ass_A(I,A)$; the completion $A^\wedge_\fp$
is still normal (\cite[Th.32.2(i)]{Mat}), and therefore its only
associated prime is $0$, so the assumption means that the radical
of $IA^\wedge_\fp$ is $\fp A^\wedge_\fp$.
Equivalently, $\dim A^\wedge_\fp/IA^\wedge_\fp=0$, so
$\dim A_\fp/IA_\fp=0$, which is the contention.
\end{example}

\begin{definition}\label{def_Lef}
Let $X$ be a noetherian scheme, $Y\subset X$ a closed subset,
$\fX$ the formal completion of $X$ along $Y$
(\cite[Ch.I, \S10.8]{EGAI}), and $f:\fX\to X$ the natural
morphism of locally ringed spaces.
We say that the pair $(X,Y)$ {\em satisfies the Lefschetz condition},
if for every open subset $U\subset X$ with $Y\subset U$, and every
locally free $\cO_U$-module $\cE$ of finite type, the natural map
$$
\Gamma(U,\cE)\to\Gamma(\fX,f^*\cE)
$$
is an isomorphism. In this case, we also say that $\Lef(X,Y)$ holds.
(Cp. \cite[Exp.X, \S2]{SGA2}.)
\end{definition}

\begin{lemma}\label{lem_was-moved}
In the situation of definition {\em\ref{def_Lef}}, suppose that
$\Lef(X,Y)$ holds, and let $U\subset X$ be any open subset such
that\/ $Y\subset U$. Then :
\begin{enumerate}
\item
The functor :
$$
\cO_{\!U}\Mod_\lfft\to\cO_\fX\Mod_\lfft
\quad :\quad
\cE\mapsto f^*\cE
$$
is fully faithful (notation of \eqref{sec_various-O-mod}).
\item
Denote by $\cO_{\!U}\Alg_\lfft$ the category of
$\cO_{\!U}$-algebras, whose underlying $\cO_{\!U}$-module
is locally free of finite type, and define likewise
$\cO_\fX\Alg_\lfft$. Then the functor :
$$
\cO_{\!U}\Alg_\lfft\to\cO_\fX\Alg_\lfft
\quad :\quad
\cA\mapsto f^*\cA
$$
is fully faithful.
\end{enumerate}
\end{lemma}
\begin{proof}(i): Let $\cE$ and $\cF$ be any two locally free
$\cO_{\!U}$-modules of finite type. We have :
$$
\Hom_{\cO_{\!U}}(\cE,\cF)=\Gamma(U,\cHom_{\cO_{\!U}}(\cE,\cF))
$$
and likewise we may compute $\Hom_{\cO_\fX}(f^*\cE,f^*\cF)$.
However, the natural map :
$$
f^*\cHom_{\cO_{\!U}}(\cE,\cF)\to\cHom_{\cO_\fX}(f^*\cE,f^*\cF)
$$
is an isomorphism of $\cO_\fX$-modules. The assertion follows.

(ii): An object of $\cO_{\!U}\Alg_\lfft$ is a locally
free $\cO_{\!U}$-module $\cA$ of finite type, together with
morphisms $\cA\otimes_{\cO_{\!U}}\cA\to\cA$ and
$1_\cA:\cO_{\!U}\to\cA$ of $\cO_{\!U}$-modules, fulfilling
the usual unitarity, commutativity and associativity conditions.
An analogous description holds for the objects of
$\cO_\fX\Alg_\lfft$, and for the morphisms of either
category. Since $\cA\otimes_{\cO_{\!U}}\cA$ is again
locally free of finite type, the assertion follows easily
from (i) : the details are left to the reader.
\end{proof}

\begin{lemma}\label{lem_Lef}
Let $A$ be a noetherian ring, $I\subset A$ an ideal,
$U\subset\Spec\,A$ an open subset, $\fU$ the formal
completion of\/ $U$ along $U\cap V(I)$.
Consider the following conditions :
\begin{enumerate}
\alphaenu
\item
$\Lef(U,U\cap V(I))$ holds.
\item
The natural map $\rho_U:\Gamma(U,\cO_{\!U})\to\Gamma(\fU,\cO_\fU)$
is an isomorphism.
\item
The natural map $\rho:A\to\Gamma(\fU,\cO_\fU)$ is an isomorphism.
\end{enumerate}
Then {\em(c)$\Rightarrow$(b)$\Leftrightarrow$(a)}, and {\em(c)} implies
that $A$ is $I$-adically complete.
\end{lemma}
\begin{proof} Clearly (a)$\Rightarrow$(b), hence we assume that
(b) holds, and we show (a). Let $V\subset U$ be an open subset
with $U\cap V(I)\subset V$ and $\cE$ a coherent locally free
$\cO_V$-module. As $A$ is noetherian, $V$ is quasi-compact,
so we may find a left exact sequence
$P_\bullet:=(0\to\cE\to\cO_V^{\oplus m}\to\cO_V^{\oplus n})$
of $\cO_V$-modules (corollary \ref{cor_ease-of-ref}).
Since the natural map of locally ringed spaces $f:\fU\to U$ is
flat, the sequence $f^*P$ is still left exact. Since the
global section functors are left exact, there follows
a ladder of left exact sequences :
$$
\Gamma(V,P_\bullet)\to\Gamma(\fU,f^*P_\bullet)
$$
which reduces the assertion to the case where $\cE=\cO_V$;
the latter is covered by the following :

\begin{claim}\label{cl_lolly}
The natural map $\rho_V:\Gamma(V,\cO_V)\to\Gamma(\fU,\cO_\fU)$
is an isomorphism.
\end{claim}
\begin{pfclaim} The isomorphism $\rho_U$ factors through $\rho_V$,
hence the latter is a surjection. Suppose that $s\in\Ker\,\rho_V$,
and $s\neq 0$; then we may find $x\in V$ such that the image $s_x$
of $s$ in $\cO_{V,x}$ does not vanish. Moreover, we may find
$a\in A$ whose image $a(x)$ in $\kappa(x)$ does not vanish, and
such that $as$ is the restriction of an element of $A$; especially,
$as\in\Gamma(U,\cO_{\!U})$, and clearly the image of $as$ in
$\Gamma(V,\cO_V)$ lies in $\Ker\,\rho_V$. Therefore, $as=0$
in $\Gamma(U,\cO_{\!U})$; however the image $as_x$ of $as$ in
$\cO_{\!U,x}$ is non-zero by construction, a contradiction.
This shows that $\rho_V$ is injective, whence the claim.
\end{pfclaim}

Finally, suppose that (c) holds; arguing as in the proof of
claim \ref{cl_lolly}, one sees that $\rho_U$ is an isomorphism.
Moreover, since $\Gamma(\fU,\cO_\fU)\simeq
\lim_{n\in\N}\,\Gamma(\fU,\cO_\fU/I^n\cO_\fU)$, the morphism
$\rho$ factors through the natural $A$-linear map $i:A\to A^\wedge$
to the $I$-adic completion of $A$. The composition with $\rho^{-1}$
yields an $A$-linear left inverse $s:A^\wedge\to A$ to $i$.
Set $N:=\Ker\,s$; clearly $s$ is surjective, hence
$A^\wedge\simeq A\oplus N$. It follows easily that $N/I^nN=0$
for every $n\in\N$, especially $N\subset\bigcap_{n\in\N}I^nA^\wedge$.
Therefore $N=0$, since $A^\wedge$ is separated for the $I$-adic
topology.
\end{proof}

\begin{proposition}\label{prop_like-Hartshorne}
Let $\phi:A\to B$ be a flat homomorphism of noetherian rings,
$I\subset A$ an ideal, $U\subset\Spec\,B$ an open subset.
Set $f:=\Spec\,\phi:\Spec\,B\to\Spec\,A$, and assume that :
\begin{enumerate}
\alphaenu
\item
$B$ is complete for the $IB$-adic topology.
\item
For every $x\in V(I)$, we have :
$\{y\in f^{-1}(x)~|~\delta(y,\cO_{\!f^{-1}(x)})=0\}\subset U$.
\item
For every $x\in\Ass_A(I,A)$, we have :
$\{y\in f^{-1}(x)~|~\delta(y,\cO_{\!f^{-1}(x)})\leq 1\}\subset U$.
\end{enumerate}
Then $\Lef(U,U\cap V(IB))$ holds.
\end{proposition}
\begin{proof} Set $\Sigma:=\Ass_A(I,A)$, and let $\fU$ be the
formal completion of $U$ along $V(IB)$; by theorem
\ref{th_brain-dead}, the $I$-preadic topology on $A$ agrees
with the $\Sigma$-symbolic topology. Let also $\cJ$ be the
family consisting of all ideals $J\subset A$ such that
$\Ass\,A/J\subset\Sigma$; it follows that the natural maps :
$$
B\to\lim_{J\in\cJ}\,B/JB
\qquad
\Gamma(\fU,\cO_\fU)\to\lim_{J\in\cJ}\,\Gamma(U,\cO_{\!U}/J\cO_{\!U})
$$
are isomorphisms (see the discussion in \eqref{subsec_dubbydubby}).
In view of lemma \ref{lem_Lef}, we are then reduced to showing :

\begin{claim} The natural map $B/JB\to\Gamma(U,\cO_{\!U}/J\cO_{\!U})$
is an isomorphism for every $J\in\cJ$.
\end{claim}
\begin{pfclaim}[] Let $f:Y:=\Spec\,B/JB\to X:=\Spec\,A/J$ be
the induced morphism; in view of corollary \ref{cor_local-depth},
it suffices to prove that $\delta(y,\cO_{\!Y})\geq 2$ whenever
$y\in Y\!\setminus\!U$. Thus, set $x:=f(y)$; by corollary
\ref{cor_depth-flat-basechange} we have :
$$
\delta(y,\cO_{\!Y})=\delta(y,\cO_{\!f^{-1}(x)})+\delta(x,\cO_{\!X}).
$$
Now, if $\delta(x,\cO_{\!X})=1$, notice that
$f(Y)\subset V(J)\subset V(I)$, by lemmata \ref{lem_ass-are-mins}(i)
and \ref{lem_little-used}(iii); hence (b) implies the contention
in this case. Lastly, if $\delta(x,\cO_{\!X})=0$, then $x\in\Ass\,A/J$
by proposition \ref{prop_misc-Ass}, hence we use assumption (c)
to conclude.
\end{pfclaim}
\end{proof}

\subsection{Cohomology of projective schemes}
\label{sec_coh-blow-up}
We begin by recalling a few generalities on graded
algebras and their homogeneous prime spectra; next
we define the blow up of a scheme along a quasi-coherent
sheaf of ideals, and we prove some basic results on
the higher direct images of quasi-coherent modules
under blow up morphisms.

\sset\subsubsection{}\label{subsec_projective-spectra}
Let $A:=\oplus_{n\in\N}A_n$ be a $\N$-graded ring, and set
$A_+:=\oplus_{n>0}A_n$, which is an ideal of $A$. Following
\cite[Ch.II, (2.3.1)]{EGAII}, one denotes by $\Proj\,A$ the
set consisting of all {\em graded prime ideals\/} of $A$
that do not contain $A_+$, and one endows $\Proj\,A$ with
the topology induced from the Zariski topology of $\Spec\,A$.
For every homogeneous element $f\in A_+$, set :
$$
D_+(f):=D(f)\cap\Proj\,A
$$
where as usual, $D(f):=\Spec\,A_f\subset\Spec\,A$. Clearly :
\set\begin{equation}\label{eq_inclusion-are-nat}
D_+(fg)=D_+(f)\cap D_+(g)
\qquad
\text{for every homogeneous $f,g\in A_+$}.
\end{equation}
The system of open subsets $D_+(f)$, for $f$ ranging over the
homogeneous elements of $A_+$, is a basis of the topology
of $\Proj\,A$ (\cite[Ch.II, Prop.2.3.4]{EGAII}), and
obviously any system of homogeneous generators
$(f_\lambda~|~\lambda\in\Lambda)$ for the ideal $A_+$ yields
an open covering
$$
\Proj\,A=\bigcup_{\lambda\in\Lambda}D_+(f_\lambda).
$$
For any homogeneous element $f\in A_+$, let also $A_{(f)}\subset A_f$
be the subring consisting of all elements of degree zero
(for the natural $\Z$-grading of $A_f$); in other words :
$$
A_{(f)}:=\sum_{k\in\N}A_k\cdot f^{-k}\subset A_f.
$$
The topological space $\Proj\,A$ carries a sheaf of rings $\cO$,
with isomorphisms of ringed spaces :
$$
\omega_f:(D_+(f),\cO_{|D_+(f)})\isom\Spec\,A_{(f)}
\qquad
\text{for every homogeneous $f\in A_+$}.
$$
and the system of isomorphisms $\omega_f$ is compatible, in
an obvious way, with the inclusions :
$$
j_{f,g}:D_+(fg)\subset D_+(f)
$$
as in \eqref{eq_inclusion-are-nat}, and with the natural
homomorphisms $A_{(f)}\to A_{(fg)}$. Especially, the locally
ringed space $(\Proj\,A,\cO)$ is a separated scheme.

\begin{example}\label{ex_classic-proj-space}
For any $d\in\N$, and any ring $R$, take $A:=R[T_0,\dots,T_d]$,
endowed with its standard $\N$-grading such that $\gr_nA$ is
the $R$-module generated by the monomials of total degree $n$,
for every $n\in\N$. The scheme $\Proj\,A$ is the
{\em projective $d$-dimensional space} over $\Spec\,R$, denoted
$$
\P^d_R.
$$
According to \eqref{subsec_projective-spectra}, it admits
the standard affine covering $\P^d_R=\bigcup_{i=0}^dD_+(T_i)$.
Set as well
$$
\tau_{ij}:=T_j/T_i
\qquad
\text{for every $i,j=0,\dots,n$}.
$$
Then clearly $A_i:=A_{(T_i)}=R[\tau_{ij}~|~j=0,\dots,d]$ is
a free polynomial $R$-algebra in $d$ variables, so that
$U_i:=D_+(T_i)$ is isomorphic to the $d$-dimensional affine
space $\A^d_R$ over $\Spec\,R$, for every $i=0,\dots,d$.
For $i\neq j$, the intersection $U_{ij}:=D_+(T_iT_j)=U_i\cap U_j$
corresponds, under these identifications, to the open
subsets
$$
\Spec\,A_i[\tau^{-1}_{ij}]\subset U_i
\qquad\text{and}\qquad
\Spec\,A_j[\tau^{-1}_{ji}]\subset U_j
$$
so we get a commutative diagram of isomorphisms of
$R$-algebras
$$
\xymatrix{ & \cO_{\P^d_R}(U_{ij}) \ar[ld] \ar[rd] \\
A_i[\tau^{-1}_{ij}] \ar[rr] & & A_j[\tau^{-1}_{ji}]
}$$
whose downward arrows are induced by the maps $\omega_{T_i}$
and $\omega_{T_j}$ of \eqref{subsec_projective-spectra}, and
where the horizontal arrow is given by the rule :
\set\begin{equation}\label{eq_change-charts}
\tau_{ik}\mapsto\tau_{jk}\cdot\tau_{ji}^{-1}
\qquad
\text{for every $k=0,\dots,d$}.
\end{equation}
\end{example}

\sset\subsubsection{}\label{subsec_in-the-situat}
Let $A':=\oplus_{n\in\N}A'_n$ be another $\N$-graded ring,
and $\phi:A\to A'$ a homomorphism of graded rings ({\em i.e.}
$\phi(A_n)\subset A'_n$ for every $n\in\N$). Following
\cite[Ch.II, (2.8.1)]{EGAII}, we let :
$$
G(\phi):=\Proj\,A'\setminus V(\phi(A_+)).
$$
This open subset of $\Proj\,A'$ is also the same as the union of all
the open subsets of the form $D_+(\phi(f))$, where $f$ ranges over
the homogeneous elements of $A_+$. The restriction to $G(\phi)$ of
$\Spec\,\phi:\Spec\,A'\to\Spec\,A$, is a continuous map
${}^a\phi:G(\phi)\to\Proj\,A$. Moreover, we have the identity :
$$
{}^a\phi^{-1}(D_+(f))=D_+(\phi(f))
\qquad\text{for every homogeneous $f\in A_+$.}
$$
Furthermore, the homomorphism $\phi$ induces a homomorphism
$\phi_{(f)}:A_{(f)}\to A'_{(\phi(f))}$, whence a morphism of schemes :
$$
\Phi_f:D_+(\phi(f))\to D_+(f).
$$
Let $g\in A_+$ be another homogeneous element; it is easily
seen that :
$$
j_{f,g}\circ\Phi_{fg}=(\Phi_f)_{|D_+(\phi(fg))}.
$$
It follows that the locally defined morphisms $\Phi_f$ glue
to a unique morphism of schemes :
$$
\Proj\,\phi:G(\phi)\to\Proj\,A
$$
such that the diagram of schemes :
\set\begin{equation}\label{eq_temperino}
{\diagram
\Spec\,A'_{(\phi(f))} \ar[d]_{\omega_{\phi(f)}}
\ar[rrr]^-{\Spec\,\phi_{(f)}} & & & \Spec\,A_{(f)} \ar[d]^{\omega_{f}} \\
D_+(\phi(f)) \ar[rrr]^{(\Proj\,\phi)_{|D_+(\phi(f))}} & & & D_+(f)
\enddiagram}
\end{equation}
commutes for every homogeneous $f\in A_+$
(\cite[Ch.II, Prop.2.8.2]{EGAII}). Lastly, we remark that,
for every homogeneous element $f'\in A'_+$, the open subset
$D_+(f')$ lies in $G(\phi)$ if and only if $f'$ lies in the
radical of the ideal of $A'$ generated by $\phi(A_+)$
(\cite[Ch.II, Cor.2.3.15]{EGAII}). Especially,
$G(\phi)=\Proj\,A'$ whenever $\phi(A_+)$ generates
the ideal $A'_+$.

\sset\subsubsection{}\label{subsec_discussion}
To ease notation, set $Y:=\Proj\,A$. Let $M:=\oplus_{n\in\Z}M_n$
be a $\Z$-graded $A$-module ({\em i.e.} $A_k\cdot M_n\subset M_{k+n}$
for every $k\in\N$ and $n\in\Z$); for every homogeneous $f\in A_+$,
denote by $M_{(f)}\subset M_f$ the submodule consisting of all
elements of degree zero (for the natural $\Z$-grading of $M_f$).
Clearly $M_{(f)}$ is an $A_{(f)}$-module in a natural way, whence
a quasi-coherent $\cO_{D_+(f)}$-module $M_{(f)}^\sim$; these
modules glue to a unique quasi-coherent $\cO_Y$-module $M^\sim$
(\cite[Ch.II, Prop.2.5.2]{EGAII}). Especially, for every $n\in\Z$,
let $A(n)$ be the $\Z$-graded $A$-module such that $A(n)_k:=A_{n+k}$
for every $k\in\Z$ (with the convention that $A_k=0$ if $k<0$). We set :
$$
\cO_Y(n):=A(n)^\sim.
$$
Any element $f\in A_n$ induces a natural isomorphism of
$D_+(f)$-modules :
$$
\cO_Y(n)_{|D_+(f)}\isom\cO_{D_+(f)}
$$
(\cite[Ch.II, Prop.2.5.7]{EGAII}). Hence, on the open subset
$$
U_n(A):=\bigcup_{f\in A_n}D_+(f)
$$
the sheaf $\cO_Y(n)$ restricts to an invertible $\cO_{\!U_n(A)}$-module.
Especially, if $A_1$ generates the ideal $A_+$, the $\cO_Y$-modules
$\cO_Y(n)$ are invertible, for every $n\in\Z$.

\begin{example}\label{ex_charts-for-O-n}
Resume the notation of example \ref{ex_classic-proj-space}.
A direct inspection of the definitions shows that $A(n)_{(T_i)}$
is the $A_i$-module generated by all fractions of the form
$T^{r_0}_0T^{r_1}_1\cdots T^{r_d}_dT_i^{-k}$, for every $k\in\N$
and every $(r_0,\dots,r_d)\in\N^{\oplus d+1}$ such that
$r_0+r_1+\cdots+r_d=n+k$. In other words, $A(n)_{(T_i)}$ is
the $A_i$-submodule of $A(n)_{T_i}$ generated by $T_i^n$.
We see then directly that this is a free $A_i$-module of
rank one. Moreover, we get a commutative diagram of
isomorphisms
$$
\xymatrix{ & \cO_{\P^d_R}(n)(U_{ij}) \ar[ld] \ar[rd] \\
T_i^nA_i[\tau^{-1}_{ij}] \ar[rr] & & T_j^nA_j[\tau^{-1}_{ji}]
}$$
whose downward arrows are induced by the restriction maps
of the sheaf $\cO_{\P^d_R}(n)$, and where the horizontal
map is the $A_i[\tau^{-1}_{ij}]$-linear map given by the rule
$$
T_i^n\mapsto T_j^n\cdot\tau^n_{ji}
$$
and where the $A_i[\tau^{-1}_{ij}]$-module structure of
$T_j^nA_j[\tau^{-1}_{ji}]$ is defined by restriction of
scalars along the isomorphism
$A_i[\tau^{-1}_{ij}]\isom A_j[\tau^{-1}_{ji}]$ as in
\eqref{eq_change-charts}.
\end{example}

\sset\subsubsection{}\label{subsec_grad_tensor}
In the situation of \eqref{subsec_in-the-situat}, let
$M:=\oplus_{n\in\Z}M_n$ be a $\Z$-graded $A$-module.
Then $M':=M\otimes_AA'$ is a $\Z$-graded $A'$-module,
with the grading defined by the rule :
$$
M'_n:=\sum_{j+k=n}\Img(M_j\otimes_\Z A'_k\to M')
\qquad\text{for every $n\in\Z$}
$$
(\cite[Ch.II, (2.1.2)]{EGAII}).
Then \cite[Ch.II, Prop.2.8.8]{EGAII} yields a natural
morphism of $\cO_{\!G(\phi)}$-modules:
$$
\nu_M:(\Proj\,\phi)^*M^\sim\to(M')^\sim_{|G(\phi)}.
$$
Moreover, set :
$$
G_1(\phi):=\bigcup_{f\in A_1}D_+(\phi(f))
$$
and notice that $G_1(\phi)\subset U_1(A')\cap G(\phi)$; by
inspecting the proof of {\em loc.cit.} we see that the restriction
$\nu_{M|G_1(\phi)}$ is an isomorphism. Especially, $\nu_M$ is an
isomorphism whenever $A_1$ generates the ideal $A_+$. It is also
easily seen that $G_1(\phi)=U_1(A')$ if $\phi(A_+)$ generates
$A'_+$.

For any $f\in A_1$, the restriction $(\nu_M)_{|D_+(\phi(f))}$
can be described explicitly : namely, we have natural
identifications
$$
\omega^*_f(M^\sim_{|D_+(f)})\isom M^\sim_{(f)}
\qquad
\omega^*_{\phi(f)}(M')^\sim_{|D_+(\phi(f))}\isom(M')^\sim_{(\phi(f))}
$$
and in view of \eqref{eq_temperino}, the morphism
$(\nu_M)_{|D_+(\phi(f))}$ is induced by the
$A'_{(\phi(f))}$-linear map :
$$
M_{(f)}\otimes_{A_{(f)}}A'_{(\phi(f))}\to M'_{(\phi(f))}
$$
given by the rule :
$$
(m_k\cdot f^{-k})\otimes(a'_j\cdot\phi(f)^{-j})\mapsto
(m_k\otimes a'_j)\cdot\phi(f)^{-j-k}
\qquad
\text{for all $k,j\in\Z$, $m_k\in M_k$, $a'_j\in A'_j$.}
$$

\sset\subsubsection{}\label{subsec_impro-somewhat}
The foregoing results can be improved somewhat, in the following
special situation. Let $R\to R'$ be a ring homomorphism,
$A$ a $\N$-graded $R$-algebra (hence the structure morphism
$R\to A$ is a ring homomorphism $R\to A_0$); the ring
$A':=R'\otimes_RA$ is naturally a $\N$-graded $R'$-algebra,
and the induced map $\phi:A\to A'$ is a homomorphism of graded
rings. In this case, obviously $\phi(A_+)$ generates the
ideal $A'_+$, hence $G(\phi)=\Proj\,A'$, and indeed, $\Proj\,\phi$
induces an isomorphism of $\Spec\,R'$-schemes :
$$
Y'\isom\Spec\,R'\times_{\Spec\,R}Y
$$
where again $Y:=\Proj\,A$ and $Y':=\Proj\,A'$. Moreover, for
every $\Z$-graded $A$-module $M$, the corresponding morphism
$\nu_M$ is an isomorphism, regardless of whether or not $A_1$
generates $A_+$ (\cite[Ch.II, Prop.2.8.10]{EGAII}). Especially,
$\nu_{\!A(n)}$ is a natural identification
(\cite[Ch.II, Cor.2.8.11]{EGAII}) :
$$
(\Proj\,\phi)^*\cO_Y(n)\isom\cO_{Y'}(n)
\qquad\text{for every $n\in\Z$}.
$$

\sset\subsubsection{}\label{subsec_Omegas}
Keep the notation of \eqref{subsec_projective-spectra},
and for every integer $d>0$, set
$$
A^{(d)}_n:=A_{nd}
\qquad\text{for every $n\in\N$, and}\qquad
A^{(d)}:=\bigoplus_{n\in\N}A^{(d)}_n.
$$
Clearly $A^{(d)}$ is an $\N$-graded ring, and a subring
of $A$, but the inclusion map $j:A^{(n)}\to A$ is not
a homomorphism of $\N$-graded rings. However, for every
homogeneous element $f\in A_+$, the map $j$ induces a
natural identification of $A$-algebras :
\set\begin{equation}\label{eq_identify-d-multiple}
A^{(d)}_{(f^d)}\isom A_{(f)}
\end{equation}
which in turn yields a natural isomorphism of $A_0$-schemes
(\cite[Ch.II, Prop.2.4.7(i)]{EGAII})
$$
\Omega^{(d)}:\Proj\,A\isom\Proj\,A^{(d)}
\qquad
\text{for every $d>0$}.
$$
To ease notation, we shall let $Y:=\Proj\,A$ and
$Y^{(d)}:=\Proj\,A^{(d)}$ for every integer $d>0$.
Likewise, if $M$ is any $\Z$-graded $A$-module, we may
consider that $\Z$-graded $A^{(d)}$-module $M^{(d)}$ whose
homogeneous direct summand of degree $n$ equals $M_{nd}$,
for every $n\in\Z$. Then the inclusion map $M^{(d)}\to M$
induces natural identifications
\set\begin{equation}\label{eq_not-quite-in-EGA}
M^{(d)}_{(f^d)}\to M_{(f)}
\qquad
\text{for every homogeneous element $f$ of $A$}.
\end{equation}
Indeed, it is easily seen that this map is surjective,
and it is injective, since it is the restriction of the
injective map $M^{(d)}_{f^d}\to M_{f^d}$. Moreover, if we
endow $M_{(f)}$ with the $A^{(d)}_{(f^d)}$-module structure
induced by \eqref{eq_identify-d-multiple}, it is easily
seen that \eqref{eq_not-quite-in-EGA} is an isomorphism
of $A^{(d)}_{(f^d)}$-modules. Thus, we obtain a natural
isomorphism of quasi-coherent $\cO_{Y^{(d)}}$-modules
$$
M^{(d)\sim}\isom\Omega^{(d)}_*M^\sim.
$$
Next, for every $n\in\Z$, let $M(n)$ be the $\Z$-graded
$A$-module such that $M(n)_k:=M_{n+k}$ for every $k\in\Z$
(with $A$-module structure deduced from that of $M$, in
the obvious way). Clearly
$$
M(nd)^{(d)}=M^{(d)}(n)
\qquad
\text{for every $n\in\Z$ and every $d>0$}
$$
whence an induced isomorphism of quasi-coherent
$\cO_{Y^{(d)}}$-modules
$$
M^{(d)}(n)^\sim\isom\Omega^{(d)}_*M(nd)^\sim
\qquad
\text{for every $n\in\Z$ and every $d>0$}.
$$
Especially, we have a natural identification
\set\begin{equation}\label{eq_twists}
\cO_{Y^{(d)}}(n)\isom\Omega^{(d)}_*\cO_Y(nd)
\qquad
\text{for every $n\in\Z$ and every $d>0$}.
\end{equation}

\sset\subsubsection{}\label{subsec_nothing-to-see}
Let $X$ be a scheme, $\cA:=\oplus_{n\in\N}\cA_n$ a $\N$-graded
quasi-coherent $\cO_{\!X}$-algebra on the Zariski site of $X$;
we let $\cA_+:=\oplus_{n>0}\cA_n$.
According to \cite[Ch.II, Prop.3.1.2]{EGAII}, there exists an
$X$-scheme $\pi:\Proj\,\cA\to X$, with natural isomorphisms
of $U$-schemes :
$$
\psi_U:U\times_X\Proj\,\cA\isom\Proj\,\cA(U)
$$
for every affine open subset $U\subset X$, and the system
of isomorphisms $\psi_U$ is compatible, in an obvious way,
with inclusions $U'\subset U$ of affine open subsets.
Especially, $\pi$ is a separated morphism.
For any integer $d>0$, every $f\in\Gamma(X,\cA_d)$ defines
an open subset $D_+(f)\subset\Proj\,\cA$, such that :
$$
D_+(f)\cap\pi^{-1}U=D_+(f_{|U})\subset\Proj\,\cA(U)
\qquad
\text{for every affine open subset $U\subset X$}.
$$

\sset\subsubsection{}\label{subsec_qcoh-mod-grad}
To ease notation, set $Y:=\Proj\,\cA$, and let again
$\pi:Y\to X$ be the natural morphism. Let $\cM:=\oplus_{n\in\Z}\cM_n$
be a $\Z$-graded $\cA$-module, quasi-coherent as a
$\cO_{\!X}$-module; for every affine open subset $U\subset X$,
the graded $\cA(U)$-module $\cM(U)$ yields a quasi-coherent
$\cO_{\!\pi^{-1}U}$-module $\cM^\sim_U$, and every inclusion of
affine open subsets $U'\subset U$ induces a natural isomorphism of
$\cO_{\!\pi^{-1}U'}$-modules : $\cM^\sim_{U|U'}\isom\cM^\sim_{U'}$.
Therefore the locally defined modules $\cM_U^\sim$ glue to
a well defined quasi-coherent $\cO_Y$-module $\cM^\sim$.

Especially, for every $n\in\Z$, denote by $\cA(n)$ the
$\Z$-graded $\cA$-module such that $\cA(n)_k:=\cA_{n+k}$ for
every $k\in\Z$, with the convention that $\cA_n=0$ whenever
$n<0$. We set:
$$
\cO_Y(n):=\cA(n)^\sim
\qquad \text{and}\qquad
\cM^\sim(n):=\cM^\sim\otimes_{\cO_Y}\cO_Y(n).
$$
Denote by $U_n(\cA)\subset Y$ the union of the open subsets
$U_n(\cA(U))$, for $U$ ranging over the affine open subsets
of $X$; from the discussion in \eqref{subsec_discussion},
it clear that the restriction $\cO_Y(n)_{|U_n(\cA)}$ is an
invertible $\cO_{\!U_n(\cA)}$-module. This open subset can
be described as follows. For every $x\in X$, let :
\set\begin{equation}\label{eq_fibres-of-A}
\cA_n(x):=\cA_{n,x}\otimes_{\cO_{\!X,x}}\kappa(x)
\qquad\text{and set}\qquad
\cA(x):=\oplus_{n\in\N}\cA_n(x)
\end{equation}
which is a $\N$-graded $\kappa(x)$-algebra; then :
\set\begin{equation}\label{eq_descr-U}
U_n(\cA)=\{y\in Y~|~\cA_n(\pi(y))\not\subset\fp(y)\}
\end{equation}
where $\fp(y)\subset\cA(\pi(y))$ denotes the prime ideal
corresponding to the point $y$.

\sset\subsubsection{}
Moreover, for every $\Z$-graded $\cA$-module $\cM$, and
every $n\in\Z$, there exists a natural morphism of
$\cO_Y$-modules :
\set\begin{equation}\label{eq_force}
\cM^\sim(n)\to\cM(n)^\sim
\end{equation}
where $\cM(n)$ is the $\Z$-graded $\cA$-module given by
the rule : $\cM(n)_k:=\cM_{n+k}$ for every $k\in\N$
(\cite[Prop.3.2.16]{EGAII}). The restriction of \eqref{eq_force}
to the open subset $U_1(\cA)$ is an isomorphism
(\cite[Ch.II, Cor.3.2.8]{EGAII}). Especially, we have
natural morphisms of $\cO_Y$-modules :
\set\begin{equation}\label{eq_annother-onne}
\cO_Y(n)\otimes_{\cO_Y}\cO_Y(m)\to\cO_Y(n+m)
\qquad
\text{for every $n,m\in\Z$}
\end{equation}
(\cite[Ch.II, Prop.3.2.6]{EGAII}) whose restrictions to
$U_1(\cA)$ are isomorphisms. Furthermore, we have a natural
morphism $\cM_0\to\pi_*\cM^\sim$ of $\cO_{\!X}$-modules
(\cite[Ch.II, (3.3.2.1)]{EGAII}), whence, by adjunction,
a morphism of $\cO_Y$-modules :
\set\begin{equation}\label{eq_force-again}
\pi^*\cM_0\to\cM^\sim.
\end{equation}
Applying \eqref{eq_force-again} to the modules $\cM_n=\cM(n)_0$,
and taking into account the isomorphism \eqref{eq_force}, we
deduce a natural morphism of $\cO_{\!U_1(\cA)}$-modules :
\set\begin{equation}\label{eq_bunny}
(\pi^*\!\cM_n)_{|U_1(\cA)}\to\cM^\sim(n)_{|U_1(\cA)}
\end{equation}
which can be described as follows. Let $U\subset X$ be
any affine open subset; for every $f\in\cA_1(U)$, the
restriction of \eqref{eq_bunny} to $D_+(f)\subset\pi^{-1}U$
is given by the morphisms
$$
\cM_n(U)\otimes_{\cO_{\!X}(U)}\cA(U)_{(f)}\to
\cM(n)(U)_{(f)}:=\sum_{k\in\Z}\cM_{k+n}(U)\cdot f^{-k}\subset\cM(U)_f.
$$
induced by the scalar multiplication
$\cM_n\otimes_{\cO_{\!X}}\cA_k\to\cM_{n+k}$. Especially,
we have natural morphisms of $\cO_Y$-modules :
\set\begin{equation}\label{eq_especiallly}
\pi^*\cA_n\to\cO_Y(n)
\qquad\text{for every $n\in\N$}
\end{equation}
whose restrictions to $U_1(\cA)$ are epimorphisms. An inspection
of the definition also shows that the diagram of $\cO_Y$-modules :
\set\begin{equation}\label{eq_diaggramma}
{\diagram \pi^*\cA_n\otimes_{\cO_Y}\pi^*\cA_m \ar[r] \ar[d] &
\pi^*\cA_{n+m} \ar[d] \\
\cO_Y(n)\otimes_{\cO_Y}\cO_Y(m) \ar[r] & \cO_Y(n+m)
\enddiagram}
\end{equation}
commutes for every $n,m\in\N$, where the top horizontal arrow
is induced by the graded multiplication
$\cA_n\otimes_{\cO_{\!X}}\cA_m\to\cA_{n+m}$, the vertical
arrows are the maps \eqref{eq_especiallly}, and the bottom
horizontal arrow is the map \eqref{eq_annother-onne}.

\sset\subsubsection{}\label{subsec_glueee}
Next, let $\cA':=\oplus_{n\in\N}\cA'_n$ be another $\N$-graded
quasi-coherent $\cO_{\!X}$-algebra on the Zariski site of $X$,
and $\phi:\cA\to\cA'$ a morphism of graded $\cO_{\!X}$-algebras;
for every affine open subset $U\subset X$, we deduce a morphism
$\phi_U:\cA(U)\to\cA'(U)$ of graded  $\cO_{\!X}(U)$-algebras,
whence an open subset $G(\phi_U)\subset\Proj\,\cA'(U)$ as
in \eqref{subsec_in-the-situat}. If $V\subset U$ is a smaller
affine open subset, the natural isomorphism
$$
V\times_U\Proj\,\cA'(U)\isom\Proj\,\cA'(V)
$$
induces an identification $V\times_UG(\phi_U)\isom G(\phi_V)$,
hence there exists a well defined open subset
$G(\phi)\subset\Proj\,\cA'$ such that the morphisms
$\Proj\,\phi_U$ glue to a unique morphism of $X$-schemes :
$$
\Proj\,\phi:G(\phi)\to\Proj\,\cA.
$$
If $\cA'_+$ is generated -- locally on $X$ -- by $\phi(\cA_+)$,
we have $G(\phi)=\Proj\,\cA'$.

Moreover, if $\cM$ is a $\Z$-graded quasi-coherent $\cA$-module,
the morphisms $\nu_{\!\cM(U)}$ assemble into a well defined morphism
of $\cO_{\!G(\phi)}$-modules :
$$
\nu_{\!\cM}:(\Proj\,\phi)^*\cM^\sim\to(\cM')^\sim_{|G(\phi)}
$$
where the grading of $\cM':=\cM\otimes_\cA\cA'$ is defined
as in \eqref{subsec_grad_tensor}. Likewise, the union of the
subsets $G_1(\phi_U)$, for $U$ ranging over the affine open
subsets of $X$, is an open subset :
\set\begin{equation}\label{eq_caveatt}
G_1(\phi)\subset U_1(\cA')\cap G(\phi)
\end{equation}
such that the restriction $\nu_{\!\cM|G_1(\phi)}$ is an isomorphism.
Especially, set $Y':=\Proj\,\cA'$; we have a natural morphism :
\set\begin{equation}\label{eq_new-nu}
\nu_{\!\cA(n)}:(\Proj\,\phi)^*\cO_Y(n)\to\cO_{Y'}(n)_{|G(\phi)}
\end{equation}
which is an isomorphism, if $\cA_1$ generates $\cA_+$
locally on $X$. Again, we have $G_1(\phi)=U_1(\cA')$ whenever
$\phi(\cA_+)$ generates $\cA'_+$, locally on $X$.

\sset\subsubsection{}\label{subsec_the-discussion}
The discussion in \eqref{subsec_impro-somewhat} implies
that any morphism of schemes $f:X'\to X$ induces a natural
isomorphism of $X'$-schemes (\cite[Ch.II, Prop.3.5.3]{EGAII}) :
\set\begin{equation}\label{eq_preimage-gralg}
\Proj\,f^*\!\cA\isom X'\times_X\Proj\,\cA
\end{equation}
and the description \eqref{eq_descr-U} implies that
\eqref{eq_preimage-gralg} restricts to an isomorphism :
$$
U_n(f^*\!\cA)\isom X'\times_XU_n(\cA)
\qquad\text{for every $n\in\N$}.
$$
Furthermore, set $Y':=\Proj\,f^*\!\cA$, and let $\pi_Y:Y'\to Y$
be the morphism deduced from \eqref{eq_preimage-gralg}; the
discussion in \eqref{subsec_impro-somewhat} implies as well that,
for any $\Z$-graded quasi-coherent $\cA$-module $\cM$, there is a
natural isomorphism :
$$
(f^*\!\cM)^\sim\isom\pi_Y^*\cM^\sim
$$
(\cite[Ch.II, Prop.3.5.3]{EGAII}). Especially, $f$ induces a
natural identification (\cite[Ch.II, Cor.3.5.4]{EGAII}) :
\set\begin{equation}\label{eq_for-arbi-n}
\cO_{Y'}(n)\isom\pi_Y^*\cO_Y(n)
\qquad\text{for every $n\in\Z$}.
\end{equation}

\sset\subsubsection{}
Keep the notation of \eqref{subsec_qcoh-mod-grad}, and
let $\cC_X$ be the category whose objects are all the pairs
$(\psi:Z\to X,\cL)$, where $\psi$ is a morphism of schemes
and $\cL$ is an invertible $\cO_{\!Z}$-module on the Zariski
site of $Z$; the morphisms $(\psi:Z\to X,\cL)\to(\psi':Z'\to X,\cL')$
are the pairs $(\beta,h)$, where $\beta:Z\to Z'$ is a morphism
of $X$-schemes, and $h:\beta^*\cL'\isom\cL$ is an isomorphism
of $\cO_{\!Z}$-modules (with composition of morphisms defined
in the obvious way). Consider the functor:
$$
F_{\!\cA}:\cC_X^o\to\Set
$$
which assigns to any object $(\psi,\cL)$ of $\cC_X$, the set
consisting of all homomorphisms of graded $\cO_{\!Z}$-algebras :
$$
g:\psi^*\!\cA\to\Sym_{\cO_{\!Z}}^\bullet\cL
$$
which are epimorphisms on the underlying $\cO_{\!Z}$-modules
(here $\Sym_{\cO_{\!Z}}^\bullet\cL$ denotes the symmetric
$\cO_{\!Z}$-algebra on the $\cO_{\!Z}$-module $\cL$);
on a morphism $(\beta,h)$ as in the foregoing, and an element
$g'\in F_{\!\cA}(\psi',\cL')$, the functor acts
by the rule :
$$
F_{\!\cA}(\beta,h)(g'):=(\Sym_{\cO_{\!Z}}^\bullet h)\circ\beta^*g'.
$$

\begin{lemma}\label{lem_repres-proj}
The object $(\pi:U_1(\cA)\to X,\cO_Y(1)_{|U_1(\cA)})$
of\/ $\cC_X$ represents the functor $F_{\!\cA}$.
\end{lemma}
\begin{proof} Given an object $(\psi:Z\to X,\cL)$
of\  $\cC_X$, and $g\in F_{\!\cA}(\psi,\cL)$, set :
$$
\P(\cL):=\Proj\,(\Sym_{\cO_{\!Z}}^\bullet\cL).
$$
According to \cite[Ch.II, Cor.3.1.7, Prop.3.1.8(iii)]{EGAII},
the natural morphism $\pi_Z:\P(\cL)\to Z$ is an isomorphism,
and clearly in this situation the natural maps
\eqref{eq_especiallly} are isomorphisms :
\set\begin{equation}\label{eq_pulllback}
\pi_Z^*\cL^{\otimes n}\isom\cO_{\P(\cL)}(n)
\qquad\text{for every $n\in\N$}.
\end{equation}
On the other hand, since $g$ is an epimorphism, we have
$G(g)=\P(\cL)$; taking \eqref{eq_preimage-gralg} into account,
we deduce a morphism of $Z$-schemes :
$$
\Proj\,g:\P(\cL)\to Y':=Z\times_X\Proj\,\cA
$$
which is the same as a morphism of $X$-schemes :
$$
\P(g):Z\to\Proj\,\cA.
$$
We need to show that the image of $\P(g)$ lies in the open
subset $U_1(\cA)$; to this aim, we may assume that both
$X$ and $Z$ are affine, say $X=\Spec\,R$, $Z=\Spec\,S$,
in which case $\cA$ is the quasi-coherent algebra associated
to a $\N$-graded $R$-algebra $A$, $\cL$ is the invertible
module associated to a projective rank one $S$-module $L$,
and $g:S\otimes_RA\to\Sym_S^\bullet L$ is a surjective homomorphism
of $R$-algebras. Then locally on $Z$, $\cL$ is generated by
elements of the form $g(1\otimes t)$, for some local sections
$t$ of $\cA_1$, and up to replacing $Z$ by an affine open
subset, we may assume that $t\in A_1$ is an element such
that $t':=g(1\otimes t)$ generates the free $S$-module $L$.
In this situation, we have $\P(\cL)=D_+(t')$, and $\P(g)$
is the same as the morphism $\Phi_t:D_+(t')\to D_+(t)$
(notation of \eqref{subsec_in-the-situat}); especially
the image of $\P(g)$ lies in $U_1(\cA)$, as required.

Moreover, the isomorphism $\pi_Z:\P(\cL)\isom Z$ is
induced by the natural identification:
\set\begin{equation}\label{eq_splitthhair}
S=S[t']_{(t')}.
\end{equation}
From this description, we also can extract an explicit
expression for $\Phi_t$; namely, it is induced by the
map of $R$-algebras :
$$
A_{(t)}\to S
\qquad\text{such that}\qquad
a_k\cdot t^{-k}\mapsto g(1\otimes a_k)\cdot t^{\prime-k}
$$
for every $k\in\N$, and every $a_k\in A_k$. Next, letting
$n:=1$ in \eqref{eq_new-nu} and \eqref{eq_for-arbi-n}, we
obtain a natural isomorphism of $\cO_{\P(\cL)}$-modules :
$$
\cO_{\P(\cL)}(1)\isom(\Proj\,g)^*\cO_{Y'}(1)\isom
(\Proj\,g)^*\circ\pi_Y^*\cO_Y(1)\isom\pi_Z^*\circ\P(g)^*\cO_Y(1)
$$
(notice that, since by assumption $g$ is an epimorphism,
we have $G_1(g)=U_1(\Sym^\bullet_{\cO_{\!Z}}\cL)=\P(\cL)$,
hence $\nu_{\!\psi^*\!\cA(1)}$ is an isomorphism).
Composition with \eqref{eq_pulllback} yields an isomorphism :
$$
h(g):\P(g)^*\cO_Y(1)\isom\cL
$$
of $\cO_{\!Z}$-modules, whence a morphism in $\cC_X$
$$
(\P(g),h(g)):(\psi,\cL)\to(\pi_{|U_1(\cA)},\cO_Y(1)_{|U_1(\cA)}).
$$
In case $X$ and $Z$ are affine, and $\cL$ is associated
to a free module $L$, generated by an element of the form
$t':=g(1\otimes t)$ as in the foregoing, we can describe
explicitly $h(g)$; namely, a direct inspection of the
construction shows that in this case $h(g)$ is induced
by the map of $S$-modules
$$
S\otimes_{A_{(t)}}A(1)_{(t)}\to L
\quad :\quad
s\otimes a_k\cdot t^{1-k}\mapsto s\cdot g(1\otimes a_k)\cdot(t')^{1-k}
\qquad\text{for every $s\in S$, $a_k\in A_k$}.
$$
Conversely, let $\beta:Z\to U_1(\cA)$ be a morphism of
$X$-schemes, and $h:\beta^*\cO_Y(1)_{|U_1(\cA)}\isom\cL$
an isomorphism of $\cO_{\!Z}$-modules. In view of the natural
isomorphisms \eqref{eq_annother-onne}, we deduce, for every
$n\in\N$, an isomorphism :
$$
h^{\otimes n}:\beta^*\cO_Y(n)_{|U_1(\cA)}\isom\cL^{\otimes n}.
$$
Combining with the epimorphisms \eqref{eq_especiallly} :
$$
\omega_n:(\pi^*\cA_n)_{|U_1(\cA)}\to\cO_Y(n)_{|U_1(\cA)}
$$
we may define the epimorphism of $\cO_{\!Z}$-modules :
\set\begin{equation}\label{eq_promise}
g(\beta,h):=
\bigoplus_{n\in\N}h^{\otimes n}\circ\beta^*(\omega_n):
\psi^*\cA\to\Sym^\bullet_{\cO_{\!Z}}\cL
\end{equation}
which, in view of \eqref{eq_diaggramma}, is a homomorphism
of graded $\cO_{\!Z}$-algebras, {\em i.e.}
$g(\beta,h)\in F(\psi,\cL)$. This homomorphism can be described
explicitly, locally on $Z$ : namely, say again that $X=\Spec\,R$,
$Z=\Spec\,S$, $\cL=L^\sim$ for a free $S$-module of rank one,
and $\cA=A^\sim$ for some $\N$-graded $R$-algebra $A$; suppose
moreover that the image of $\beta$ lies in an open subset
$D_+(t)\subset U_1(\cA)$, for some $t\in A_1$. Then $\beta$
comes from a ring homomorphism $\beta^\natural:A_{(t)}\to S$,
$h$ is an $S$-linear isomorphism $S\otimes_{A_{(t)}}A(1)_{(t)}\isom L$,
and $t':=h(1\otimes t)$ is a generator of $L$; moreover, $\omega_n$
is the epimorphism deduced from the map :
$$
A_n\otimes_RA_{(t)}\to A(n)_{(t)}
\quad :\quad
a_n\otimes b_k\cdot t^{-k}\mapsto a_nb_k\cdot t^{-k}
\qquad\text{for every $a_n\in A_n$, $b_k\in A_k$}
$$
By inspecting the construction, we see therefore that $g$
is the direct sum of the morphisms :
$$
g_n:S\otimes_RA_n\to L^{\otimes n}
\quad :\quad
s\otimes a_n\mapsto
s\cdot\beta^\natural(a_n\cdot t^{-n})\cdot t^{\prime\otimes n}
\qquad
\text{for every $s\in S$, $a_n\in A_n$}.
$$
Finally, it is easily seen that the natural transformations :
\set\begin{equation}\label{eq_rules}
g\mapsto(\P(g),h(g))
\qquad\text{and}\qquad
(\beta,h)\mapsto g(\beta,h)
\end{equation}
are inverse to each other : indeed, the verification can be
made locally on $Z$, hence we may assume that $X$ and $Z$ are
affine, and $\cL$ is free, in which case one may use the explicit
formulae provided above.
\end{proof}

\begin{definition}\label{def_blow-up-well}
Let $X$ be a scheme, $\cI$ a quasi-coherent sheaf of
ideals of $\cO_{\!X}$.
\begin{enumerate}
\item
The {\em Rees algebra} of $\cI$ is the quasi-coherent
$\N$-graded $\cO_{\!X}$-algebra
$$
\sR(\cI)_\bullet:=\bigoplus_{n\in\N}\cI^n
$$
with multiplication law deduced in the obvious way from
that of $\cO_{\!X}$.
\item
The {\em blowing up} of $\cI$ is the $X$-scheme
$$
\Proj\,\sR(\cI)_\bullet\to X.
$$
\end{enumerate}
\end{definition}

\begin{remark}\label{rem_blowing-up}
(i)\ \
With the notation of definition \ref{def_blow-up-well},
let $\pi:Y\to X$ be the blowing up of $\cI$. Notice that
$\cI$ is an invertible $\cO_{\!X}$-module if and only if
it is generated, locally on $X$, by a regular section
of $\cO_{\!X}$, and then the natural map
$$
\Sym^n_{\cO_X}\cI\to\cI^n
$$
is an isomorphism, for every $n\in\N$. Taking into
account  \cite[Ch.II, Cor.3.1.7, Prop.3.1.8(iii)]{EGAII},
it follows that if $\cI$ is an invertible ideal, then the
blowing up of $\cI$ is an isomorphism.

(ii)\ \
Let $X'\to X$ be any flat morphism of schemes; in light
of \eqref{subsec_the-discussion}, we get a natural
isomorphism from $X'\times_XY\to X'$ to the blowing up
of the quasi-coherent ideal $\cI\cO_{\!X'}$ of $\cO_{\!X'}$.

(iii)\ \
Especially, let $U\subset X$ be the complement in $X$
of the closed subscheme $\Spec\,\cO_{\!X}/\cI$. From
(i) and (ii) we deduce that the restriction $\pi^{-1}U\to U$
of $\pi$ is an isomorphism.

(iv)\ \
Notice that $\sR(\cI)_\bullet$ is generated, locally on
$X$, by its degree one direct summand $\cI$, hence
$\cO_Y(1)$ is an invertible $\cO_Y$-module. On the
other hand, a simple inspection shows that
$$
\sR(\cI)_\bullet(n)=\cI^n\cdot\sR(\cI)_\bullet
\qquad
\text{for every $n\in\N$}
$$
whence a natural identification
$$
\cO_Y(n)\isom\cI^n\cdot\cO_Y
\qquad
\text{for every $n\in\N$}.
$$
\end{remark}

\begin{proposition}\label{prop_univ-prop-blow-up}
With the notation of definition {\em\ref{def_blow-up-well}},
the blowing up $Y\to X$ of $\cI$ is characterized, up to
unique isomorphism of $X$-schemes, by the following two
conditions :
\begin{enumerate}
\item
The sheaf of ideals $\cI\cdot\cO_Y$ is an invertible
$\cO_Y$-module.
\item
For every $X$-scheme $Z$ such that $\cI\cdot\cO_{\!Z}$
is an invertible $\cO_{\!Z}$-module, there exists a
unique morphism $Z\to Y$ of $X$-schemes.
\end{enumerate}
\end{proposition}
\begin{proof} The uniqueness up to unique isomorphism
of an $X$-scheme $\pi:Y\to X$ fulfilling conditions (i) and
(ii) is clear. It remains to check that these conditions
hold for $Y:=\Proj\,\sR(\cI)_\bullet$.

However, condition (i) follows immediately from remark
\ref{rem_blowing-up}(iv). Next, let $f:Z\to X$ be as in
(ii), and set $\cL:=\cI\cdot\cO_{\!Z}$; the induced map
$f^*\cO_{\!X}\to\cO_{\!Z}$ induces an epimorphism of
$\cO_{\!Z}$-modules
$$
\phi_n:f^*\cI^n\to\cI^n\cO_{\!Z}\isom\Sym^n_{\cO_Z}\cL
\qquad
\text{for every $n\in\N$}
$$
and it is clear that the system $(\phi_n~|~n\in\N)$ amounts
to a morphism of $\N$-graded $\cO_{\!Z}$-algebras
$$
\phi_\bullet:f^*\sR(\cI)_\bullet\to\Sym^\bullet_{\cO_Z}\cL
$$
which in turn corresponds to a morphism $g:Z\to Y$ of
$X$-schemes, by virtue of lemma \ref{lem_repres-proj}.
Lastly, let $g':Z\to Y$ be any other morphism of
$X$-schemes, and denote by $U\subset X$ (resp.
$U'\subset Z$) the complement in $X$ (resp. in $Z$)
of the closed subset $\Spec\,\cO_{\!X}/\cI$ (resp.
$\Spec\,\cO_{\!Z}/\cI\cO_{\!Z}$); clearly $f(U')\subset U$,
hence $g$ and $g'$ both restrict to morphisms
$U'\to\pi^{-1}U$ of $X$-schemes. Then remark
\ref{rem_blowing-up}(iii) implies that
$$
g_{|U'}=g'_{|U'}.
$$
Next, since $\cI\cO_{\!Z}$ is generated, locally on $Z$,
by a regular section of $\cO_{\!Z}$, the open subset $U'$
is {\em schematically dense} in $Z$
(\cite[Ch.IV, D\'ef.11.10.2]{EGAIV-3}), and on the other
hand, we know that $\pi$ is a separated morphism (see
\eqref{subsec_nothing-to-see}); in view of
\cite[Ch.IV, Prop.11.10.1]{EGAIV-3} we deduce that
$g=g'$, which concludes the verification of (ii).
\end{proof}

\sset\subsubsection{}
We may generalize as follows remark \ref{rem_blowing-up}(ii).
Let $X\to X_0$ be a morphism of schemes, $\cI_0\subset\cO_{\!X_0}$
a quasi-coherent ideal, and set $\cI:=\cI_0\cO_{\!X}$; let
also $\pi_0:Y_0\to X_0$ and $\pi:Y\to X$ be the blowing up
morphisms of $\cI_0$ and resepctively $\cI$. By proposition
\ref{prop_univ-prop-blow-up} there exists a unique morphism
of schemes $\psi:Y\to Y_0$ that makes commute the diagram
\set\begin{equation}\label{eq_diagr-of-blowups}
{\diagram Y \ar[r]^-\psi \ar[d]_\pi &
Y_0 \ar[d]^{\pi_0} \\
X \ar[r] & X_0.
\enddiagram}
\end{equation}
Then we have :

\begin{lemma}\label{lem_blow-up-as-cl-immers}
The morphism of $X$-schemes $\bar\psi:Y\to Y':=Y_0\times_{X_0}X$
induced by $\psi$ is a closed immersion, and identifies
$Y$ with $\Spec\,\cO_{Y'}/\cJ$, where
$$
\cJ:=\bigcup_{n\in\N}\Ann_{\cO_{Y'}}(\cI^n).
$$
\end{lemma}
\begin{proof} Let $Z\to X$ be any morphism of schemes
such that $\cI\cO_{\!Z}$ is an invertible $\cO_{\!Z}$-module;
since $\cI=\cI_0\cO_{\!X}$, proposition
\ref{prop_univ-prop-blow-up} yields a unique morphism
of $X_0$-schemes $Z\to Y_0$, so there exists a unique
morphism of $X$-schemes $g:Z\to Y'$ as well. Next, since
$\cI\cO_{\!Z}$ is invertible, it is easily seen that the
induced morphism $\cO_{\!Y'}\to g_*\cO_{\!Z}$ factors
through $\cO_{Y'}/\cJ$, whence a unique morphism of
$X$-schemes $Z\to Y'':=\Spec\,\cO_{Y'}/\cJ$. Lastly,
by construction $\cI\cO_{Y'}$ is a locally principal
ideal of $\cO_{\!Y'}$, and then it is clear that
$\cI\cO_{Y''}$ is an invertible $\cO_{Y''}$-module.
Summing up, we see that the morphism $Y''\to X$
enjoys the universal property that characterizes
the blowing up of $\cI$, whence the lemma.
\end{proof}

\sset\subsubsection{}\label{subsec_cohomol-of-blow-up}
Let $A$ be a ring, $r\geq 1$ an integer, $\bff:=(f_1,\dots,f_r)$
a sequence of elements of $A$, and $I\subset A$ the ideal
generated by $\bff$. We view $A$ as an algebra over the ring
$A_0:=\Z[T_1,\dots,T_r]$, via the ring homomorphism
$$
g:A_0\to A
\qquad
T_i\mapsto f_i
\qquad
i=1,\dots,r
$$
and denote by $I_0\subset A_0$ the ideal generated by
$(T_1,\dots,T_r)$. Set also $X:=\Spec\,A$, $X_0:=\Spec\,A_0$,
let $\cI\subset\cO_{\!X}$ (resp. $\cI_0\subset\cO_{\!X_0}$) be
the quasi-coherent sheaf of ideals arising from $I$ (resp.
from $I_0$), and $\pi:Y\to X$ (resp. $\pi_0:Y_0\to X_0$) the
blowing up of $\cI$ (resp. of $\cI_0$). Clearly
$\cI_0\cO_Y=\cI$, whence a commutative diagram
\eqref{eq_diagr-of-blowups}, with bottom horizontal
arrow given by $\Spec\,g:X\to X_0$. Furthermore, notice
that we have systems of monomorphisms of invertible
$\cO_Y$-modules :
\set\begin{equation}\label{eq_inverse-sys-blowup}
\cO_{Y_0}(n+1)\to\cO_{Y_0}(n)
\qquad
\cO_Y(n+1)\to\cO_Y(n)
\qquad
\text{for every $n\in\N$}
\end{equation}
corresponding to the inclusion maps $\cI^{n+1}\cO_{Y_0}\to\cI^n\cO_{Y_0}$
and $\cI^{n+1}\cO_Y\to\cI^n\cO_Y$ under the natural identifications of
remark \ref{rem_blowing-up}(iv).

\begin{proposition}\label{prop_compl-secant-blow-up}
With the notation of \eqref{subsec_cohomol-of-blow-up},
suppose moreover that the sequence $\bff$ is completely
secant. Then, for every $n\in\N$ we have :
\begin{enumerate}
\item
The natural map $I^n\to H^0(Y,\cO_Y(n))$ is an isomorphism.
\item
$H^p(Y,\cO_Y(n))=0$ for every $p>0$.
\end{enumerate}
\end{proposition}
\begin{proof} We consider first the ring $A_0$ and the sequence
$\bff_0:=(T_1,\dots,T_r)$; we endow $A_0$ with its standard
$\N$-grading, such that $\gr_nA_0$ is the $\Z$-module generated
by the monomials of total degree $n$, for every $n\in\N$. Let
also $\sR_0\subset A_0[V]$ be the Rees algebra associated with
the $I_0$-adic filtration of $A_0$, so that $Y_0=\Proj\,\sR_0$
(see definition \ref{def_Rees}(iii); but since here we
have a descending filtration on $A_0$, we let $V:=U^{-1}$,
so $\gr_n\sR_0:=V^nI_0^n$ for every $n\in\N$). We consider
the morphism of $\N$-graded $\Z$-algebras
$$
h:A_0\to\sR_0
\qquad
T_i\mapsto VT_i
\qquad
i=1,\dots,r. 
$$
Notice that $\phi(A_{0+})$ generates the ideal $\sR_{0+}$
of $\sR_0$, hence $\phi:=\Proj\,h$ is well defined on the
whole of $Y_0$ (see \eqref{subsec_in-the-situat}), and
we get a commutative diagram of morphisms of schemes
$$
\xymatrix{
Y_0 \ar[r]^-\phi \ar[d]_{\pi_0} & \P^{r-1}_\Z \ar[d] \\
X_0 \ar[r] & \Spec\,\Z.
}$$
For every $i=1,\dots,r$, set $U_i:=D_+(T_i)\subset\P^{r-1}_\Z$,
and $U'_i:=\phi^{-1}U_i$. Recall that there are natural
isomorphisms $\omega_i:
\cO_{\P_\Z^{r-1}}(U_i)\isom A_i:=\Z[\tau_{ij}~|~j=1,\dots,r]$
(notation of example \ref{ex_classic-proj-space}); we
wish to give a corresponding description of
the $A_i$-algebra $B_i:=\cO_{Y_0}(U'_i)$. To this aim,
notice that, by definition, every element of
$B_i\subset\Z[T_1,\dots,T_r,V,(VT_i)^{-1}]$ is a finite
sum of terms of the form
\set\begin{equation}\label{eq_terms}
(VT_1)^{\alpha_1}\cdots(VT_r)^{\alpha_r}\cdot(VT_i)^{-k}
\cdot P(T_1,\dots,T_r)=\tau_{i1}^{\alpha_1}\cdots\tau_{ir}^{\alpha_r}
\cdot P(T_i\tau_{i1},\dots,T_i\tau_{ir})
\end{equation}
where $k\in\N$ is any integer, $P\in A$ any polynomial, and
$\alpha_1+\cdots+\alpha_r=k$; {\em i.e.}
$$
B_i=A_i[T_i]
\qquad
\text{for every $i=1,\dots,r$}.
$$
Likewise, for every $n\in\N$, the elements of $\sR_0(n)_{(T_i)}$
are the finite sums of terms \eqref{eq_terms} where $k\in\N$
is any integer, $P\in A$ any polynomial, and
$\alpha_1+\cdots+\alpha_r=n+k$; {\em i.e.}
$$
\sR_0(n)_{(T_i)}=(VT_i)^nB_i
\qquad
\text{for every $i=1,\dots,r$}.
$$
For every $i\neq j$, set also $U_{ij}:=U_i\cap U_j$ and
$U'_{ij}:=\phi^{-1}U_{ij}$, and recall that the isomorphisms
$\omega_i$ induce isomorphisms
$\omega_{ij}:\cO_{\P_\Z^{r-1}}(U_{ij})\isom A_i[\tau^{-1}_{ij}]$,
such that the composition
$$
\psi_{ij}:=\omega_{ji}\circ\omega^{-1}_{ij}:A_i[\tau^{-1}_{ij}]\isom
A_j[\tau^{-1}_{ji}]
$$
is given by the rule \eqref{eq_change-charts}. There follow
natural identifications $\cO_{Y_0}(U'_{ij})\isom B_i[\tau^{-1}_{ij}]$,
whence a commutative diagram
$$
\xymatrix{ A_i[\tau^{-1}_{ij}] \ar[r]^-{\psi_{ij}} \ar[d] &
A_j[\tau^{-1}_{ji}] \ar[d] \\
B_i[\tau_{ij}^{-1}] \ar[r]^-{\psi'_{ij}} & B_j[\tau^{-1}_{ji}]
}$$
whose vertical arrows are induced by the map
$\cO_{\P^{r-1}_\Z}\to\phi_*\cO_{Y_0}$ associated with $\phi$,
and where $\psi'_{ij}$ is the isomorphism given by the rule :
$$
T_i\mapsto T_j\cdot\tau_{ji}.
$$
Likewise, we get a commutative diagram
$$
\xymatrix{ & \cO_{Y_0}(n)(U'_{ij}) \ar[ld] \ar[rd] \\
(VT_i)^nB_i[\tau^{-1}_{ij}] \ar[rr] & & (VT_j)^nB_j[\tau^{-1}_{ji}]
}$$
whose downward arrows are induced by the restriction maps
of the sheaf $\cO_{Y_0}(n)$, and where the horizontal
map is the $B_i[\tau^{-1}_{ij}]$-linear map given by the rule
$$
(VT_i)^n\mapsto(VT_j)^n\cdot\tau^n_{ji}
$$
and where the $B_i[\tau^{-1}_{ij}]$-module structure of
$(VT_j)^nB_j[\tau^{-1}_{ji}]$ is defined by restriction of
scalars along the isomorphism $\psi'_{ij}$. Comparing
with example \ref{ex_charts-for-O-n} we deduce a natural
isomorphism of $\cO_{\P^{r-1}_\Z}$-modules :
$$
\phi_*\cO_{Y_0}(n)\isom\bigoplus_{k\in\N}\cO_{\P^{r-1}_\Z}(n+k)
\qquad
\text{for every $n\in\N$}.
$$
On the other hand, from \cite[Ch.III, Prop.2.1.12]{EGAIII}
we get isomorphisms of $\N$-graded $\Z$-modules:
$$
H^p\Bigl(\P^{r-1}_\Z,\bigoplus_{k\in\N}\cO_{\P^{r-1}_\Z}(k)\Bigr)=
         \left\{\begin{array}{lll}
           0 & \qquad & \text{for every $p>0$} \\
           A & \qquad & \text{for $p=0$}
         \end{array}\right.
$$
Now, since $\phi$ is affine, we have
$$
H^p(Y_0,\cO_{Y_0}(n))=H^p(\P^{r-1}_\Z,\phi_*\cO_{Y_0}(n))
\qquad
\text{for every $n\in\N$}
$$
and assertion (ii) already follows. We also deduce an isomorphism
of $\N$-graded $\Z$-modules
$$
H^0(Y_0,\cO_{Y_0}(n))\isom
H^0\Bigl(\P^{r-1}_\Z,\bigoplus_{k\in\N}\cO_{\P^{r-1}_\Z}(n+k)\Bigr)=I_0^n
$$
where $\gr_kI_0^n:=\gr_{n+k}A$, for every $k\in\N$.
To conclude the proof in this case, it remains only to
check that this identification is the inverse of the
natural map of (i). To this aim, we consider the induced
commutative diagram
$$
\xymatrix{ H^0(Y_0,\cO_{Y_0}(n)) \ar[r] \ar[d] &
H^0(\P^{r-1}_\Z,\phi_*\cO_{Y_0}(n)) \ar[d] \\
H^0(U'_i,\cO_{Y_0}(n)) \ar[r] & H^0(U_i,\phi_*\cO_{Y_0}(n)).
}$$
Now, let $T^\alpha:=T_1^{\alpha_1}\cdots T_r^{\alpha_r}\in\gr_kI^n$
be any monomial (for any $k\in\N$); by inspecting the
foregoing notation, we see that the image of $T^\alpha$
in $H^0(U'_i,\cO_{Y_0}(n))$ equals
$(VT_i)^nT^k_i\tau_{i,1}^{\alpha_1}\cdots\tau_{i,r}^{\alpha_r}$,
which maps to the section
$T_i^{n+k}\tau_{i,1}^{\alpha_1}\cdots\tau_{i,r}^{\alpha_r}$ of
$H^0(U_i,\phi_*\cO_{Y_0}(n))$. But the latter is also the
image of $T^\alpha$ under the restriction map
$H^0(\P^{r-1}_\Z,\cO_{\P^{r-1}}(n+k))\to H^0(U_i,\cO_{\P^{r-1}}(n+k))$,
whence the contention.

Next, let $A$, $\bff$ and $I$ as in
\eqref{subsec_cohomol-of-blow-up}, and we define likewise
$\sR\subset A[V]$ as the Rees algebra associated with the
$I$-adic filtration of $A$; we consider the homomorphism
of graded $\Z$-algebras
$$
G:\sR_0\to\sR
\qquad
V^nx\mapsto V^ng(x)
\qquad
\text{for every $n\in\N$ and every $x\in I^n_0$}.
$$
Notice that $G$ is even a morphism of $\N$-graded
$A_0$-algebras, for the $A_0$-algebra structure on
$A$ and $\sR$ induced by $g$. There follows a morphism
of $\N$-graded $A$-algebras
$$
h':A\otimes_{A_0}\sR_0\to\sR
$$
whose restriction $\gr_nh'$ is induced by the natural
map $A\otimes_{A_0}I_0^n\to I^n$, for every $n\in\N$.
The latter is an isomorphism for every $n\in\N$, by
virtue of lemma \ref{lem_secant-badabum}, and it is
clear that the morphism $\psi$ of
\eqref{subsec_cohomol-of-blow-up} equals $\Proj\,h'$;
in view of \eqref{subsec_impro-somewhat}, we deduce
that the diagram \eqref{eq_diagr-of-blowups} is
cartesian under the current assumptions, and moreover
$h'$ induces natural isomorphisms of $\cO_Y$-modules
$$
\psi^*\cO_{Y_0}(n)\isom\cO_Y(n)
\qquad
\text{for every $n\in\N$}.
$$
Denote by $\fU'$ the affine open covering $(U'_i~|~i=1,\dots,r)$
of $Y_0$, and set $\psi^{-1}\fU':=(\psi^{-1}U'_i~|~i=1,\dots,r)$.
there follows a natural isomorphism of alternating \v{C}ech
complexes
\set\begin{equation}\label{eq_pull-back-Cech}
\bar C{}^\bullet_\alt(\fU',\cO_{Y_0}(n))\otimes_{A_0}A\isom
\bar C{}^\bullet_\alt(\psi^{-1}\fU',\cO_Y(n))
\qquad
\text{for every $n\in\N$}.
\end{equation}

\begin{claim}\label{cl_use-first-E_1}
The natural map
$\bar C{}^\bullet_\alt(\fU',\cO_{Y_0}(n))\derotimes_{A_0}A\to
\bar C{}^\bullet_\alt(\fU',\cO_{Y_0}(n))\otimes_{A_0}A$
is an isomorphism in $\sD(A\Mod)$, for every $n\in\N$.
\end{claim}
\begin{pfclaim} We consider the standard spectral sequence
$$
E_1^{pq}:=\Tor_{-q}^{A_0}(\bar C{}^p_\alt(\fU',\cO_{Y_0}(n)),A)
\Rightarrow
H^{p+q}(\bar C{}^\bullet_\alt(\fU',\cO_{Y_0}(n))\derotimes_{A_0}A)
$$
and notice that it suffices to check that $E_1^{pq}=0$ for
every $p\neq 0$, since in that case the abutment is naturally
isomorphic to the cohomology of the complex
$(E_1^{0,\bullet},d_1^{0,\bullet})$, which is the same as
$\bar C{}^\bullet_\alt(\fU',\cO_{Y_0}(n))\otimes_{A_0}A$.
However, with the foregoing notation, we see that
$\bar C{}^q_\alt(\fU',\cO_{Y_0}(n))$ is a direct sum of
finitely many $A_0$-modules, each of which is a localization
of $B_i$, for some $i\leq r$, so we are reduced to checking
that $\Tor_p^{A_0}(B_i,A)=0$ for $i=1,\dots,r$ and every $p>0$.
But in turn, $B_i$ is the degree $0$ summand of a $\Z$-graded
$A_0$-module, which is a localization of $\sR_0$, so we are
further reduced to checking that $\Tor_p^{A_0}(\sR_0,A)=0$
for every $p>0$, or equivalently, that $\Tor_p^{A_0}(I_0^n,A)=0$
for every $p>0$ and every $n\in\N$. The latter follows from
lemma \ref{lem_secant-badabum}.
\end{pfclaim}

Combining \eqref{eq_pull-back-Cech}, claim
\ref{cl_use-first-E_1} and theorem \ref{th_Cech-resolve}(ii),
we get a natural isomorphism :
$$
\bar C{}^\bullet_\alt(\fU',\cO_{Y_0}(n))\derotimes_{A_0}A\isom
R\Gamma(Y,\cO_Y(n))
\qquad
\text{in $\sD(A\Mod)$, for every $n\in\N$}.
$$
By the same token, we have a natural isomorphism
$\bar C{}^\bullet_\alt(\fU',\cO_{Y_0}(n))\isom
R\Gamma(Y_0,\cO_{Y_0}(n))$ in $\sD(A_0\Mod)$, whence
a standard spectral sequence for every $n\in\N$
$$
E(n)^{pq}_2:=\Tor^{A_0}_{-p}(H^q(Y_0,\cO_{Y_0}(n)),A)\Rightarrow
H^{p+q}(\bar C{}^\bullet_\alt(\fU',\cO_{Y_0}(n))\derotimes_{A_0}A)
$$
and the previous case tells us that $E(n)^{pq}_2=0$ for
every $q>0$ and every $n\in\N$, and moreover
$E(n)^{p0}_2=\Tor^{A_0}_{-p}(I^n_0,A)$ for every $p\in\Z$ and
every $n\in\N$. Lastly, by invoking again lemma
\ref{lem_secant-badabum} we see that $E(n)^{p0}_2=0$
whenever $p\neq 0$, and $E(n)^{00}_2=I^n_0A=I^n$ for every
$n\in\N$, which concludes the proof of the proposition.
\end{proof}

\begin{theorem}\label{th_too-difficult-for-me}
With the notation of \eqref{subsec_cohomol-of-blow-up},
the following conditions are equivalent :
\begin{enumerate}
\alphaenu
\item
The ring $A$ satisfies condition $\mathrm{(a)^{un}_\bff}$.
\item
The inverse system $(H^p(Y,\cO_Y(n))~|~n\in\N)$ deduced from
the system of maps \eqref{eq_inverse-sys-blowup} is uniformly
essentially zero for every $p>0$, and the same holds for the
kernel and cokernel of the natural morphism of inverse systems
$(I^n\to H^0(Y,\cO_Y(n))~|~n\in\N)$.
\item
There exists $n\in\N$ such that $I^n\cdot H^p(Y,\cO_Y)=0$ for
every $p>0$, and $I^n$ annihilates the kernel and cokernel
of the natural map $A\to H^0(Y,\cO_Y)$.
\end{enumerate}
\end{theorem}
\begin{proof} Define $I_0\subset A$, the Rees algebras
$\sR_0$, $\sR$, and the morphism $\pi_0:Y_0\to X_0:=\Spec\,A_0$
as in the proof of proposition \ref{prop_compl-secant-blow-up}.
We set $Y':=Y_0\times_{X_0}X$ and define $\cJ\subset\cO_{Y'}$,
$\psi:Y\to Y$ and $\bar\psi:Y\to Y'$ as in lemma
\ref{lem_blow-up-as-cl-immers}, so that $\cJ$ is the kernel
of the induced map $\bar\psi{}^\sharp:\cO_{Y'}\to\bar\psi_*\cO_Y$.
We shall consider as well the condition :

\begin{itemize}
\item[(d)]
There exists $k\in\N$ such that
$I^k\!\cJ=I^k\cdot\cTor^{X_0}_i(\cO_{Y_0},\cO_{\!X})=0$
for every $i>0$.
\end{itemize}

\begin{claim}\label{cl_a-implies-local-tor}
If condition (a) holds, there exists $k\in\N$ such that
$$
I^k\cdot\Tor^{A_0}_i(I_0^n,A)=I^k\cdot\Ker(I_0^n\otimes_{A_0}A\to I^n)=0
\qquad
\text{for every $n\in\N$ and every $i>0$}.
$$
\end{claim}
\begin{pfclaim} For every integer $p>0$, let $k(p)$ be the
step of the uniformly essentially zero system
$(\Tor_p^{A_0}(A_0/I_0^n,A)~|~n\in\N)$. For $p>1$ (resp. for
$p=1$), the integer $k(p)$ is also the step of the uniformly
essentially zero system
$(T^n_{p-1}:=\Tor_{p-1}^{A_0}(I_0^n,A)~|~n\in\N)$
(resp. the uniformly essentially zero system
$(K^n:=\Ker(I_0^n\otimes_{A_0}A\to I^n)~|~n\in\N)$). Now,
for every $a\in I^{k(p)}_0$, scalar multiplication by $a$
on $I^n_0$ factors through the inclusion map
$I^{n+k(p)}_0\to I_0^n$, so it induces the zero map on
$T^n_{p-1}$ (resp. on $K^n$) for every $n\in\N$. Lastly,
recall that the homological dimension of $A_0$ equals
$r+1$, so $T^n_p=0$ for every $p>r+1$. We may therefore
choose $k:=\max(k(1),\dots,k(r+1))$.
\end{pfclaim}

(a)$\Rightarrow$(d): Indeed, recall that we have a natural
$\cO_{Y'}(U\times_{X_0}X)$-linear isomorphism
$$
\Gamma(U\times_{X_0}X,\cTor^{X_0}_i(\cO_{Y_0},\cO_{\!X}))
\isom
T_i(U):=\Tor^{A_0}_i(H^0(\cO_{Y_0}(U),A)
$$
for every affine open subsets $U\subset Y_0$.
Since $Y'$ is quasi-compact, and taking into account lemma
\ref{lem_blow-up-as-cl-immers}, it will then suffice to find
$k\in\N$ and a finite open covering
$U_\bullet:=(U_\lambda~|~\lambda\in\Lambda)$ of $Y_0$ such that
\begin{itemize}
\item
$I^k\cdot T_i(U_\lambda)=0$ for every $i>0$ and every
$\lambda\in\Lambda$
\item
$I^k\cdot\Ker(\bar\psi{}^\sharp_{U_\lambda}:
\cO_{Y_0}(U_\lambda)\otimes_{A_0}A\to\cO_Y(\bar\psi{}^{-1}U_\lambda))=0$
for every $\lambda\in\Lambda$.
\end{itemize}
However, we may choose $U_\bullet$ such $\cO_{Y_0}(U_\lambda)$
is a direct summand of a graded localization of $\sR_0$,
for every $\lambda\in\Lambda$; then $T_i(U_\lambda)$ is a
direct summand of a localization of $\Tor_i^{A_0}(\sR_0,A)$,
and $\bar\psi{}^\sharp_{U_\lambda}$ is a direct summand of a
localization of the natural map $\sR_0\otimes_{A_0}A\to\sR$.
Thus, the assertion follows from claim
\ref{cl_a-implies-local-tor}.

(b)$\Rightarrow$(c): The natural identification
$\cI^n\cO_Y\isom\cO_Y(n)$ of remark \ref{rem_blowing-up}(iv)
yields exact cohomology sequences
$$
H^p(Y,\cO_Y(n))\xrightarrow{\ \alpha^p_n\ }H^p(Y,\cO_Y)\to
H^p(Y,\cO_Y/\cI^n\cO_Y)
\qquad
\text{for every $p,n\in\N$}.
$$
However, assumption (b) implies that for every $p>0$
there exists $k\in\N$ such that $\alpha^p_n$ is the
zero map, for every $n\geq k$; hence, for such $k$
the term $H^p(Y,\cO_Y)$ is a submodule of
$H^p(Y,\cO_Y/\cI^k\cO_Y)$, which is obviously annihilated
by $I^k$. By the same token we get a commutative ladder
with exact rows for every $n\in\N$ :
$$
\xymatrix{ 0 \ar[r] & I^n \ar[r] \ar[d]_{\beta_n} &
A \ar[r] \ar[d]^{\gamma_n} & A/I^n \ar[r] \ar[d]^{\delta_n} & 0 \\
0 \ar[r] & H^0(Y,\cO_Y(n)) \ar[r] & H^0(Y,\cO_Y) \ar[r] &
H^0(Y,\cO_Y/\cI^n\cO_Y)
}$$
whence an exact sequence
$$
0\to\Ker\,\beta_n\xrightarrow{\ \tau_n\ }\Ker\,\gamma_n\to
\Ker\,\delta_n\to\Coker\,\beta_n\xrightarrow{\ \sigma_n\ }
\Coker\,\gamma_n\to\Coker\,\delta_n.
$$
However, assumption (b) implies that there exists $n\in\N$
such that $\tau_n$ and $\sigma_n$ both vanish; since
$\Ker\,\delta_n$ and $\Coker\,\delta_n$ are both annihilated
by $I^n$, we conclude that (c) holds.

Next, we show that (d) implies both (a) and (b). To this
aim, we remark :

\begin{claim}\label{cl_inverse-Tor}
Suppose that condition (d) holds; then we have :
\begin{enumerate}
\item
The inverse system
$(\cT^i_n:=\cTor^{X_0}_i(\cO_{Y_0}(n),\cO_{\!X})~|~n\in\N)$
induced by the system of morphisms \eqref{eq_inverse-sys-blowup}
is uniformly essentially zero, for every $i>0$.
\item
Likewise, the inverse system
$(\cJ(n):=\Ker(\cO_{Y'}(n)\to\bar\psi_*\cO_Y(n))~|~n\in\N)$
is uniformly essentially zero.
\end{enumerate}
\end{claim}
\begin{pfclaim}(i): There exists an affine open covering
$(U_1,\dots,U_r)$ of $Y_0$ such that
$$
(\cI\cO_{Y_0})_{|U_i}=f_i\cO_{Y_0}
\qquad
\text{for every $i=1,\dots,r$}
$$
and $f_i$ is a regular element of $\cO_{Y_0}(U_i)$, whence
an identification
$$
(\cI^n\cO_{Y_0})_{|U_i}=f^n_i\cO_{\!U_i}\isom\cO_{\!U_i}
\qquad
\text{for every $i=1,\dots,r$}.
$$
Under this identification, the restriction
$\cO_{Y_0}(n+1)_{|U_i}\to\cO_{Y_0}(n)_{|U_i}$ of the morphism
\eqref{eq_inverse-sys-blowup} corresponds to the scalar
multiplication by $f_i$. Hence, the inverse system
$((\cT^i_n)_{|U_i}~|~n\in\N)$ is naturally identified with
the inverse system $(\cT'^i_n~|~n\in\N)$ with
$\cT'^i_n:=\cTor_i^{X_0}(\cO_{U_i},\cO_{\!X})$ for every
$i,n\in\N$, with transition maps given by scalar
multiplication by $f_i$, whence (i).

For the proof of (ii) one argues in the same way, using
the fact that $I^k\!\cJ=0$ for some $k\in\N$ whose existence
is ensured by condition (d) : the details shall be left
to the reader.
\end{pfclaim}

We consider now the inverse systems of spectral sequences
provided by proposition \ref{prop_local-Tors}
\set\begin{equation}\label{eq_sink}
E(n)_2^{pq}:=H^p(Y_0\times_{X_0}X,\cTor_{-q}^{X_0}(\cO_{Y_0}(n),\cO_{\!X}))
\Rightarrow H^{p+q}(A\derotimes_{A_0}R\Gamma\cO_{Y_0}(n))
\end{equation}
where, as usual, the transition maps $E(n+1)_2^{pq}\to E(n)_2^{pq}$
are induced by the morphisms \eqref{eq_inverse-sys-blowup}.
Clearly $E_2^{pq}=0$ if either $p<0$ or $q>0$, and claim
\ref{cl_inverse-Tor}(i) implies that the inverse system
$(E(n)_2^{pq}~|~n\in\N)$ is uniformly essentially zero
whenever $q<0$. It follows already that the inverse system
$(E(n)_\infty^{pq}~|~n\in\N)$ is uniformly essentially zero
whenever $p+q<0$. On the other hand, proposition
\ref{prop_compl-secant-blow-up} shows that
$$
H^{-i}(A\derotimes_{A_0}R\Gamma\cO_{Y_0}(n))=\Tor_i^{A_0}(A,I^n_0)
\qquad
\text{for every $i\in\Z$}.
$$
Summing up, and taking into account lemma
\ref{lem_inverse-Serre-subcat}, a simple induction shows
that the system $(\Tor_i^{A_0}(A,I^n_0)~|~n\in\N)$ is uniformly
essentially zero for every $i>0$. Notice also that
$E(n)_2^{p0}=H^p(Y',\cO_{Y'}(n)$ for every $n,p\in\N$.

\begin{claim}\label{cl_again-usual}
(i)\ \
The inverse system $(E(n)_s^{p0})~|~n\in\N)$ is uniformly
essentially zero, for every $p>0$ and every $s\geq 2$.

(ii)\ \
The kernel and cokernel of the natural morphism
of inverse systems
$$
h_\bullet:(A\otimes_{A_0}I_0^n~|~n\in\N)\to(E(n)^{00}_2~|~n\in\N)
$$
are both uniformly essentially zero.
\end{claim}
\begin{pfclaim}(i): We argue by descending induction on $p$ :
if $p\geq r$, we have $E(n)_2^{p0}=0$ for every $n\in\N$,
by theorem \ref{th_Cech-resolve}(iii), hence the claim
trivially holds in this case. Thus, suppose that $q<r$,
and that the claim is already known for every $p>q$ and
every $s\geq 2$. We show, by descending induction on
$s\geq 2$, that $(E(n)_s^{q0}~|~n\in\N)$ is uniformly
essentially zero. Since $E(n)_2^{q-s,s-1}=0$ for every
$s\geq 2$, we have $E(n)^{q-s,s-1}_s=E(n)_s^{q+s,1-s}=0$ as
well, whenever $q+s\geq r$, and therefore
$E(n)_s^{q0}=E(n)_\infty^{q0}=0$ for every $p>0$ and every
$s\geq r$. Thus, suppose that $2\leq t<r$, and that the
assertion is known for every $s>t$; we consider the exact
sequence of inverse systems
$$
0\to (E(n)^{q0}_{t+1}~|~n\in\N)\to(E(n)_t^{q0}~|~n\in\N)\to
(E(n)^{q+t,1-t}_t~|~n\in\N).
$$
Since $1-t<0$, the third term is uniformly essentially
zero, and the same holds for the first term, by inductive
assumption; by lemma \ref{lem_inverse-Serre-subcat}, we
deduce that the same holds for the middle term, as required.

(ii): The map $h_n$ is the composition of the
projection $h'_n:A\otimes_{A_0}I^n_0\to E(n)^{00}_\infty$ with
the injective map $h''_n:E(n)^{00}_\infty\to E(n)^{00}_2$ deduced
from the spectral sequences \eqref{eq_sink}. But since the
inverse system $(E(n)_2^{p,-p}~|~n\in\N)$ is uniformly
essentially zero for every $p>0$, and vanishes for $p\geq r$,
the same holds for the inverse systems $(E(n)_\infty^{p,-p}~|~n\in\N)$,
and consequently the kernel of $h'_\bullet$ is uniformly essentially
zero. Likewise, it is easily seen that the cokernel of
$h''_\bullet$ is uniformly essentially zero.
\end{pfclaim}

Define $\cJ(n)$ as in claim \ref{cl_inverse-Tor}(ii), and
set $J(n)^p:=H^p(Y',\cJ(n)$ for every $n,p\in\N$; we get
an exact sequence of inverse systems
$$
(J(n)^p~|~n\in\N)\to(E(n)_2^{p0}~|~n\in\N)
\xrightarrow{\ k^p_\bullet\ }(H^p(Y,\cO_Y(n))~|~n\in\N)
\to(J(n)^{p+1}~|~n\in\N)
$$
for every $p\in\N$, and claims \ref{cl_inverse-Tor}(ii) and
\ref{cl_again-usual}(i) imply that the first and third terms
are uniformly essentially zero, whenever $p>0$. Therefore,
the same holds for the middle term. Lastly, we consider the
commutative diagram
$$
\xymatrix{ (A\otimes_{A_0}I^n_0~|~n\in\N) \ar[r]^-{m_\bullet}
\ar[d]_{h_\bullet} & (I^n~|~n\in\N) \ar[d]^{l_\bullet} \\
(E(n)_2^{00}~|~n\in\N) \ar[r]^-{k^0_\bullet} &
(H^0(Y,\cO_Y(n))~|~n\in\N)
}$$
where $m_\bullet$ is given by the natural surjections, and
$h_\bullet$ is defined as in claim \ref{cl_again-usual}.
Since the inverse system $(J(n)^p~|~n\in\N)$ is uniformly
essentially zero for $p=0,1$, we see that the kernel and
cokernel of $k^0_\bullet$ are both uniformly essentially zero.
In view of claim \ref{cl_again-usual}(ii), it follows
immediately that the kernel of $m_\bullet$ is uniformly
essentially zero, and the same holds for the kernel and
cokernel of $l_\bullet$, which concludes the verification
of both (a) and (b).

Lastly, we check that (c) implies both (a) and (b). To this
aim we consider, for every $n\in\N$, the spectral sequences
provided by proposition \ref{prop_local-Tors}(ii)
$$
\begin{aligned}
E(n)_2^{pq}:=\,& \Tor_{-p}^{A_0}(I_0^n,H^q(Y,\cO_Y))\Rightarrow
H^{p+q}(I_0^n\derotimes_{A_0}R\Gamma\cO_Y) \\
F(n)_2^{pq}:=\,& H^p(Y,\cTor^{X_0}_{-q}(\cI_0^n,\cO_Y))
\Rightarrow H^{p+q}(I_0^n\derotimes_{A_0}R\Gamma\cO_Y).
\end{aligned}
$$
Moreover, for every $\in\N$ denote by
$$
\phi_n:\cTor_0^{X_0}(\cO_{Y_0},\cI_0^n)=\pi_0^*\cI_0^n\to
\cI^n_0\cO_{Y_0}=\cO_{Y_0}(n)
$$
the natural morphism, and notice that
$\psi^*\pi_0^*\cI_0^n=\cTor^{X_0}_0(\cI_0^n,\cO_Y)$.

\begin{claim}\label{cl_religion}
(i)\ \
The inverse system $(\Ker\,\phi_n~|~n\in\N)$ is
uniformly essentially zero, and the same holds for
the inverse system $(\cTor^{X_0}_q(\cI_0^n,\cO_Y)~|~n\in\N)$,
for every $q>0$.

(ii)\ \
The inverse system
$(H^i(I_0^n\derotimes_{A_0}R\Gamma\cO_Y)~|~n\in\N)$ is
uniformly essentially zero for $i\neq 0$.

(iii)\ \
The spectral sequences $E(n)_2^{\bullet\bullet}$ induce,
for every $p\in\Z$, an epimorphism of inverse systems
$(E(n)^{p0}_2~|~n\in\N)\to(E(n)^{p0}_\infty~|~n\in\N)$
whose kernel is uniformly essentially zero.

(iv)\ \
The spectral sequences $F(n)_2^{\bullet\bullet}$ induce,
for every $p\in\Z$, a monomorphism of inverse systems
$(F(n)^{p0}_\infty~|~n\in\N)\to(F(n)^{p0}_2~|~n\in\N)$
whose cokernel is uniformly essentially zero.

(v)\ \
The inverse systems $(E(n)_2^{p,-p}~|~n\in\N)$ and
$(F(n)_2^{-p,p}~|~n\in\N)$ are uniformly essentially zero
for $p<0$, and vanish identically for $p>0$.
\end{claim}
\begin{pfclaim}(i): We look at the spectral sequence
$$
G(n)^2_{pq}:=\cTor^{Y_0}_p(\cO_Y,\cTor_q^{X_0}(\cO_{Y_0},\cI_0^n))
\Rightarrow\cTor^{X_0}_{p+q}(\cI_0^n,\cO_Y)
$$
of proposition \ref{prop_local-Tors}(i). Since $Y_0$ is a
noetherian scheme, corollary \ref{cor_a-unif-for-noether}(iii)
implies that the inverse system
$(\cTor_q^{X_0}(\cO_{Y_0},\cI_0^n)~|~n\in\N)$ is uniformly
essentially zero for every $q>0$, so the same holds for
the inverse system $(G(n)_2^{pq}~|~n\in\N)$, for every $p\in\N$
and every $q>0$. By the same token, the inverse system
$(\Ker\,\phi_n~|~n\in\N)$ is uniformly essentially zero;
since the functor
$$
\cTor^{Y_0}_p(\cO_Y,-):\sD^-(\cO_{Y_0}\Mod_\qcoh)\to
\sD^-(\cO_Y\Mod_\qcoh)
$$
is triangulated for every $p\in\Z$, we get a short exact
sequence
$$
0\to\cTor^{Y_0}_p(\cO_Y,\Ker\,\phi_n)\to
G(n)^2_{p0}\to\cTor^{Y_0}_p(\cO_Y,\cO_{Y_0}(n))\to 0
\qquad
\text{for every $p,n\in\N$}
$$
which shows that the inverse system $(G(n)^2_{p0}~|~n\in\N)$
is uniformly essentially zero for every $p>0$. The
assertion now follows by a simple induction, using lemma
\ref{lem_inverse-Serre-subcat}.

(ii): Condition (c) implies that the inverse system
$(E(n)_2^{pq}~|~n\in\N)$ is uniformly essentially zero
for $q>0$, and clearly it vanishes identically when
either $q<0$ or $p>0$. It also vanishes identically
for $q\geq r$, due to theorem \ref{th_Cech-resolve}(iii).
Summing up, we see that the filtration induced by the
spectral sequence $E(n)^{\bullet\bullet}_2$ on its abutment
$H^i:=H^i(I_0^n\derotimes_{A_0}R\Gamma\cO_Y)$ is finite for
every $i\in\Z$, and the inverse systems of its graded
subquotients $(E(n)_\infty^{pq}~|~n\in\N)$ are uniformly
essentially zero whenever $p+q>0$, whence the assertion
in case $i>0$, after invoking as usual lemma
\ref{lem_inverse-Serre-subcat}. For $i<0$, we argue
likewise, using the spectral sequences
$F(n)_2^{\bullet\bullet}$ : indeed, we see again easily
that the filtration induced on the abutment $H^i$ is
finite for every $i\in\Z$; moreover, if $p<0$ the
inverse system $(F(n)_2^{pq}~|~n\in\N)$ is identically
zero, and it is uniformly essentially zero if $q<0$,
by virtue of (i), whence the assertion.

(iii): For $p>0$ we have $E(n)_2^{p0}=0$ for every
$n\in\N$, so the assertion is trivial in this case.
For $p\leq 0$, notice that $E(n)_s^{p+s,1-s}=0$ for every
$s\geq 2$, hence $E(n)_s^{p0}$ is a quotient of $E(n)_2^{p0}$
for every $s\geq 2$, whence a surjection
$E(n)_2^{p0}\to E(n)_\infty^{p0}$ as sought. Moreover,
notice that the inverse system $(E(n)_2^{p-s,s-1}~|~n\in\N)$
is uniformly essentially zero for every $s\geq 2$, due
to condition (c) and lemma \ref{lem_ess-zero-Tors}, and
it vanishes identically for $s\geq r$. Then the same holds
for the inverse system $(E(n)_s^{p-s,s-1}~|~n\in\N)$, for
every $s\geq 2$, and the assertion follows easily.

(iv): For $p<0$ we have $F(n)_2^{p0}=0$ for every
$n\in\N$, so the assertion is trivial in this case.
For $p\geq 0$, notice that $F(n)_s^{p-s,s-1}=0$ for
every $s\geq 2$, hence $F(n)_s^{p0}$ is a submodule
of $F(n)_2^{p0}$ for every $s\geq 2$, whence a injection
$F(n)_\infty^{p0}\to F(n)_2^{p0}$ as sought. Moreover,
notice that the inverse system $(F(n)_2^{p+s,1-s}~|~n\in\N)$
is uniformly essentially zero for every $s\geq 2$, due
to (i), and it vanishes identically for $s\geq r$. Then
the same holds for the inverse system
$(F(n)_s^{p+s,1-s}~|~n\in\N)$, for every $s\geq 2$, and
the assertion follows easily.

(v): The assertion for $p>0$ is obvious, and the case
$p<0$ for $(E(n)_2^{p,-p}~|~n\in\N)$ follows from condition
(c) and lemma \ref{lem_ess-zero-Tors}, and for
$(F(n)_2^{-p,p}~|~n\in\N)$ follows from (i).
\end{pfclaim}

Assumption (c) implies that the natural map $A\to H^0(Y,\cO_Y)$
induces a morphism of inverse systems
$$
(\Tor^{A_0}_{-p}(I^n_0,A)~|~n\in\N)\to(E(n)_2^{p0}~|~n\in\N)
\qquad
\text{for every $p\in\Z$}
$$
whose kernel and cokernel are uniformly essentially
zero. Combining with claim \ref{cl_religion}(ii,iii),
it follows already that the inverse system
$(\Tor^{A_0}_p(A/I^n_0,A)~|~n\in\N)$ is uniformly
essentially zero for every $p>1$.
Next, since $\cO_{Y_0}(n)$ is a flat $\cO_{Y_0}$-module
for every $n\in\N$, we have an exact sequence of
$\cO_Y$-modules
$$
0\to\psi^*\Ker\,\phi_n\to\cTor^{X_0}_0(\cI_0^n,\cO_Y)\to
\psi^*\cO_{Y_0}(n)=\cO_Y(n)\to 0
$$
whence exact sequences
$$
H^p(Y,\psi^*\Ker\,\phi_n)\to F(n)_2^{p0}\to
H^p(Y,\cO_Y(n))\to H^{p+1}(Y,\psi^*\Ker\,\phi_n)
$$
for every $p\in\Z$ and every $n\in\N$. Combining with
claim \ref{cl_religion}(ii,iv), we deduce that the inverse
system $(H^p(Y,\cO_Y(n))~|~n\in\N)$ is uniformly essentially
zero for every $p>0$.

Lastly, claim \ref{cl_religion}(v) implies that
$(E(n)_\infty^{p,-p}~|~n\in\N)$ and
$(F(n)_\infty^{-p,p}~|~n\in\N)$ vanish identically
for every $p>0$, and are uniformly essentially zero
inverse systems for $p<0$. It follows that the spectral
sequences $E(n)_2^{\bullet\bullet}$ and $F(n)_2^{\bullet\bullet}$
induce respectively a monomorphism and an epimorphism
$$
(E(n)_\infty^{00}~|~n\in\N)\to
(H^0(I_0^n\derotimes_{A_0}R\Gamma\cO_Y)~|~n\in\N)
\to(F(n)^{00}_\infty~|~n\in\N)
$$
and the kernel and cokernel of the composition of these
morphisms are uniformly essentially zero inverse systems.
Taking into account claim \ref{cl_religion}(iv) we get
therefore a commutative diagram of inverse systems
$$
\xymatrix{ (E(n)_2^{00}~|~n\!\in\!\N) \ar[r] \ddouble &
(E(n)_\infty^{00}~|~n\in\N) \ar[r] &
(F(n)^{00}_\infty~|~n\in\N) \ar[d] \\
(I_0^n\otimes_{A_0}H^0(Y,\cO_Y)~|~n\in\N) \ar[r]^-{\beta_\bullet}
&(H^0(Y,\cO_Y(n)~|~n\in\N) & (F(n)^{00}_2~|~n\in\N) \ar[l]
}$$
all whose arrows, except possibly $\beta_\bullet$, have
uniformly essentially zero kernels and cokernels, but
then also $\beta_\bullet$ does.

\begin{claim}\label{cl_its-the-right-map}
$\beta_\bullet$ is the map induced by the identification
$\omega_n:\cI^n\cO_Y\isom\cO_Y(n)$ of remark
\ref{rem_blowing-up}(iv) and the natural map
$\psi^*\pi^*_0\cI^n_0\to\cI^n\cO_Y$.
\end{claim}
\begin{pfclaim} Indeed, by construction, the map
$F(n)^{00}_2\to H^0(Y,\cO_Y(n))$ is obtained from
the identification $\omega_n$ and the natural
isomorphism $\cTor_0^{X_0}(\cO_{Y_0},\cI^n_0)\isom
\psi^*\pi^*_0\cI^n_0$, so we are reduced to checking
that the map
$$
E(n)^{00}_2=I^n_0\otimes_{A_0}H^0(Y,\cO_Y)\to
F(n)^{00}_2=H^0(Y,\psi^*\pi^*_0\cI^n_0)
$$
deduced from the foregoing commutative diagram, agrees
with the standard one that assigns to any element
$x\otimes s$ its class in $H^0(Y,\psi^*\pi^*_0\cI^n_0)$.
To this aim, set $U_i:=D_+(f_i)\subset Y$ for $i=1,\dots,r$,
and let $Y'':=\amalg_{i=1}^rU_i$, the affine $X$-scheme
given by the disjoint union of the affine open subsets
$U_1,\dots,U_r$. We have an obvious morphism $\nu:Y''\to Y$
of schemes, that restricts to the inclusion map on each
open subset $U_i$ of $Y''$, and let also
$\psi'':=\psi\circ\nu:Y''\to Y_0$. We have as well two
systems of spectral sequences $E''(n)_2^{\bullet\bullet}$
and $F''(n)_2^{\bullet\bullet}$, obtained as in the foregoing,
after replacing $Y$ by $Y''$ and $\cO_Y$ by $\cO_{Y''}$.
By functoriality of the spectral sequences, there follows
a commutative diagram
$$
\xymatrix{
E(n)^{00}_2=I^n_0\otimes_{A_0}H^0(Y,\cO_Y) \ar[r] \ar[d] &
F(n)^{00}_2=H^0(Y,\psi^*\pi^*_0\cI^n_0) \ar[d] \\
E''(n)^{00}_2=I^n_0\otimes_{A_0}H^0(Y'',\cO_{Y''}) \ar[r] &
F''(n)^{00}_2=H^0(Y'',\psi''^*\pi^*_0\cI^n_0)
}$$
and notice that the right vertical arrow is injective.
Thus, it suffices to check that the bottom horizontal
arrow is the expected map. However, we have
$E''(n)_2^{00}=E''(n)_\infty^{00}=E''(n)^0_\infty=F''(n)^0_\infty=
F''(n)^{00}_\infty=F''(n)^{00}_2$ for every $n\in\N$, and a
direct inspection shows that the bottom horizontal map is
the identity map, whence the claim.
\end{pfclaim}

From claim \ref{cl_its-the-right-map} and condition (c),
it follows already that the natural morphism of inverse
systems
$(I_0^n\otimes_{A_0}A~|~n\in\N)\to(H^0(Y,\cO_Y(n))~|~n\in\N)$
has uniformly essentially zero kernel and cokernel.
But the latter factors as the composition of the epimorphism
of inverse systems $(I_0^n\otimes_{A_0}A~|~n\in\N)\to(I^n~|~n\in\N)$
and the natural morphism
$(I^n~|~n\in\N)\to(H^0(Y,\cO_Y(n))~|~n\in\N)$. It follows that
the kernels and cokernels of both of the latter morphisms
are uniformly essentially zero, and the proof of the theorem
is concluded.
\end{proof}

\begin{remark}\label{rem_too-difficult-for-me}
By direct inspection of the proof, we obtain the following
refinement of theorem \ref{th_too-difficult-for-me}. For
every $r\in\N$ there exists an integer $\nu(r)$ such that,
if the ring $A$ satisfies condition $\mathrm{(a)^{un}_\bff}$
with step $\leq r$, then the inverse system
$(H^p(Y,\cO_Y(n))~|~n\in\N)$ is uniformly essentially
zero with step $\leq\nu(r)$, and the same holds for
the kernel and cokernel of the inverse system
$(I^n\to H^0(Y,\cO_Y(n))~|~n\in\N)$. The detailed
verification shall be left to the reader.
\end{remark}

\section{Duality theory}\label{chap_duality}

\subsection{Duality for quasi-coherent modules}
Let $f:X\to Y$ be a morphism of schemes; theorem
\ref{th_trivial-dual} falls short of proving that $Lf^*$
is a left adjoint to $Rf_*$, since the former functor
is defined on bounded above complexes, while the latter
is defined on complexes that are bounded below.
This deficency has been overcome by N.Spaltenstein's
paper \cite{Sp}, where he shows how to extend the usual
constructions of derived functors to unbounded complexes.
On the other hand, in many cases one can also construct
a {\em right adjoint\/} to $Rf_*$. This is the subject
of Grothendieck's duality theory. We shall collect the
statements that we need from that theory, and refer
to the original source \cite{Har} for the rather elaborate
proofs.

\sset\subsubsection{}\label{subsec_Groth-dual}
For any morphism of schemes $f:X\to Y$, let
$\bar f:(X,\cO_{\!X})\to(Y,f_*\cO_{\!X})$ be the corresponding
morphism of ringed spaces. The functor
$$
\bar f_*:\cO_{\!X}\Mod\to f_*\cO_{\!X}\Mod
$$
admits a left adjoint $\bar f{}^*$ defined as usual
(\cite[Ch.0, \S4.3.1]{EGAI}). In case $f$ is affine, $\bar f$
is flat and the unit and counit of the adjunction restrict to
isomorphisms of functors :
$$
\begin{aligned}
\one\to\bar f_*\circ\bar f{}^* &
:f_*\cO_{\!X}\Mod_\qcoh\to f_*\cO_{\!X}\Mod_\qcoh \\
\bar f{}^*\circ\bar f_*\to\one &
:\cO_{\!X}\Mod_\qcoh\to\cO_{\!X}\Mod_\qcoh.
\end{aligned}
$$
on the corresponding subcategories of quasi-coherent modules
(\cite[Ch.II, Prop.1.4.3]{EGAII}). It follows easily that, on
these subcategories, $\bar f{}^*$ is also a right adjoint to
$\bar f_*$.

\begin{lemma}\label{lem_Groth-dual}
Keep the notation of \eqref{subsec_Groth-dual}, and suppose
that $f$ is affine. Then the functor :
$$
\bar f{}^*:\sD^+(f_*\cO_{\!X}\Mod)_\qcoh\to\sD^+(\cO_{\!X}\Mod)_\qcoh
$$
is right adjoint to the functor
$R\bar f_*:\sD^+(\cO_{\!X}\Mod)_\qcoh\to\sD^+(f_*\cO_{\!X}\Mod)_\qcoh$.
Moreover, the unit and counit of the resulting adjunction
are isomorphisms of functors.
\end{lemma}
\begin{proof} By trivial duality (theorem \ref{th_trivial-dual}),
$\bar f{}^*$ is left adjoint to $R\bar f_*$ on $\sD^+(\cO_{\!X}\Mod)$;
it suffices to show that the unit (resp. and counit) of this latter
adjunction are isomorphisms for every object $K^\bullet$ of
$\sD^+(\cO_{\!X}\Mod)_\qcoh$ (resp. $L^\bullet$ of
$\sD^+(f_*\cO_{\!X}\Mod)_\qcoh$).
Concerning the unit, we can use a Cartan-Eilenberg injective
resolution of $\bar f{}^*K^\bullet$, and a standard spectral
sequence, to reduce to the case where $K^\bullet=\cF[0]$ for
a quasi-coherent $f_*\cO_{\!X}$-module $\cF$; however, the natural
map $\bar f_*\circ\bar f{}^*\cF\to R\bar f_*\circ\bar f{}^*\cF$
is an isomorphism, so the assertion follows from
\eqref{subsec_Groth-dual}.
\end{proof}

\begin{remark}
(i)\ \ Even when both $X$ and $Y$ are affine, neither the unit
nor the counit of adjunction in \eqref{subsec_Groth-dual} is an
isomorphism on the categories of all modules. For a counterexample
concerning the unit, consider the case of a finite injection of
domains $A\to B$, where $A$ is local and $B$ is semi-local (but
not local); let $f:X:=\Spec\,B\to\Spec\,A$
be the corresponding morphism, and denote by $\cF$ the
$f_*\cO_{\!X}$-module supported on the closed point, with stalk
equal to $B$. Then one verifies easily that
$\bar f_*\circ\bar f{}^*\cF$ is a direct product of finitely
many copies of $\cF$, indexed by the closed points of $X$.
Concerning the counit, keep the same morphism $f$, let
$U:=X\!\setminus\!\{x\}$, where $x$ is a closed point and set
$\cF:=j_!\cO_{\!U}$; then it is clear that $\bar f{}^*\circ\bar f_*\cF$
is supported on the complement of the closed fibre of $f$.

(ii)\ \ Similarly, one sees easily that $\bar f{}^*$ is not
a right adjoint to $\bar f_*$ on the category of all
$f_*\cO_{\!X}$-modules.

(iii)\ \ Incidentally, the analogous functor $\bar f{}^*_{\!\et}$
defined on the {\'e}tale (rather than Zariski) site is a right
adjoint to $\bar f_{\!\et *}:\cO_{\!X,\et}\Mod\to f_*\cO_{\!X,\et}\Mod$,
and on these sites the unit and counit of the adjunctions are
isomorphisms for all modules.
\end{remark}

\sset\subsubsection{}
For the construction of the right adjoint $f^!$ to $Rf_*$,
one considers first the case where $f$ factors as a composition
of a finite morphism followed by a smooth one, in which case
$f^!$ admits a corresponding decomposition. Namely :
\begin{itemize}
\item
In case $f$ is smooth, we shall consider the functor :
$$
f^\sharp:\sD(\cO_Y\Mod)\to\sD(\cO_{\!X}\Mod)\qquad
K^\bullet\mapsto f^*K^\bullet\otimes_{\cO_{\!X}}
\Lambda^n_{\cO_{\!X}}\Omega^1_{X/Y}[n]
$$
where $n$ is the locally constant relative dimension function of $f$.
\item
In case $f$ is finite, we shall consider the functor :
$$
f^\flat:\sD^+(\cO_Y\Mod)\to\sD^+(\cO_{\!X}\Mod)
\qquad K^\bullet\mapsto
\bar f{}^*R\cHom^\bullet_{\cO_Y}(f_*\cO_{\!X},K^\bullet).
$$
\end{itemize}
If $f$ is any quasi-projective morphism, such factorization
can always be found locally on $X$, and one is left with the
problem of patching a family of locally defined functors
$((f_{|U_i})^!~|~i\in I)$, corresponding
to an open covering $X=\bigcup_{i\in I} U_i$. Since such
patching must be carried out in the derived category, one
has to take care of many cumbersome complications.

\sset\subsubsection{}
Recall that a finite morphism $f:X\to Y$ is said to be
{\em pseudo-coherent\/} if $f_*\cO_{\!X}$ is a pseudo-coherent
$\cO_Y$-module. This condition is equivalent to the pseudo-coherence
of $f$ in the sense of \cite[Exp.III, D{\'e}f.1.2]{SGA6}.

\begin{lemma}\label{lem_rightist}
Let $f:X\to Y$ be a finite morphism of schemes. Then :
\begin{enumerate}
\item
If $f$ is finitely presented, the functor
$$
\cF\mapsto\bar f{}^*\cHom_{\cO_Y}(f_*\cO_{\!X},\cF)
$$
is right adjoint to $f_*:\cO_{\!X}\Mod_\qcoh\to\cO_Y\Mod_\qcoh$.
\item
If $f$ is pseudo-coherent, the functor
$$
K^\bullet\mapsto\bar f{}^*R\cHom^\bullet_{\cO_Y}(f_*\cO_{\!X},K^\bullet)
\quad :\quad
\sD^+(\cO_Y\Mod)_\qcoh\to\sD^+(\cO_{\!X}\Mod)_\qcoh
$$
is right adjoint to
$$
Rf_*:\sD^+(\cO_{\!X}\Mod)_\qcoh\to\sD^+(\cO_Y\Mod)_\qcoh
$$
(notation of \eqref{sec_various-O-mod}).
\end{enumerate}
\end{lemma}
\begin{proof} (i): Under the stated assumptions, $f_*\cO_{\!X}$
is a finitely presented $\cO_Y$-module (by claim \ref{cl_fin-present}),
hence the functor
$$
\cF\mapsto\cHom_{\cO_Y}(f_*\cO_{\!X},\cF)
\quad :\quad \cO_Y\Mod\to f_*\cO_{\!X}\Mod
$$
preserves the subcategories of quasi-coherent modules, hence
it restricts to a right adjoint for the forgetful functor
$f_*\cO_{\!X}\Mod_\qcoh\to\cO_Y\Mod_\qcoh$ (claim
\ref{cl_right-to-forget}). Then the assertion follows from
\eqref{subsec_Groth-dual}.

(ii) is analogous : since $f_*\cO_{\!X}$ is pseudo-coherent,
the functor
$$
K^\bullet\mapsto R\cHom_{\cO_Y}^\bullet(f_*\cO_{\!X},K^\bullet)
\quad :\quad \sD^+(\cO_Y\Mod)\to\sD^+(f_*\cO_{\!X}\Mod)
$$
preserves the subcategories of complexes with quasi-coherent
homology, hence it restricts to a right adjoint for the forgetful
functor $\sD^+(f_*\cO_{\!X}\Mod)_\qcoh\to\sD^+(\cO_Y\Mod)_\qcoh$, by
lemma \ref{lem_triv-dual}(iii). To conclude, it then suffices
to apply lemma \ref{lem_Groth-dual}.
\end{proof}

\begin{proposition}\label{prop_sharp-flat}
Suppose that
$X\xrightarrow{\ f\ }Y\xrightarrow{\ g\ }Z\xrightarrow{\ h\ }W$
are morphisms of schemes. Then :
\begin{enumerate}
\item
If $f$ and $g$ are finite and $f$ is pseudo-coherent, there
is a natural isomorphism of functors on $\sD^+(\cO_{\!Z}\Mod)$ :
$$
\xi_{f,g}:(g\circ f)^\flat\isom f^\flat\circ g^\flat.
$$
\item
If $f$ and $g\circ f$ are finite and pseudo-coherent, and
$g$ is smooth of bounded fibre dimension, there is a natural
isomorphism of functors on $\sD^+(\cO_{\!Z}\Mod)$ :
$$
\zeta_{f,g}: (g\circ f)^\flat\isom f^\flat\circ g^\sharp.
$$
\item
If $f$, $g$ and $h\circ g$ are finite and pseudo-coherent,
and $h$ is smooth of bounded fibre dimension, then the diagram
of functors on $\sD^+(\cO_W\Mod)$ :
$$
\xymatrix{
(h\circ g\circ f)^\flat \ar[rr]^-{\zeta_{g\circ f,h}}
\ar[d]_{\xi_{f,h\circ g}} & &
(g\circ f)^\flat\circ h^\sharp \ar[d]^{\xi_{f,g}\circ h^\sharp} \\
f^\flat\circ(h\circ g)^\flat \ar[rr]^-{f^\flat(\zeta_{g,h})} & &
f^\flat\circ g^\flat\circ h^\sharp}
$$
commutes.
\item
If, moreover
$$
\xymatrix{
Y \ar[r]^-g \ar[d]_{h'} & Z \ar[d]^-h \\
W' \ar[r]^-j & W
}$$
is a cartesian diagram of schemes, such that $h$ is smooth
of bounded fibre dimension and $j$ is finite and pseudo-coherent,
then there is a natural isomorphism of functors :
$$
\vartheta_{j,h}:g^\flat\circ h^\sharp\isom h^{\prime\sharp}\circ
j^\flat.
$$
\end{enumerate}
\end{proposition}
\begin{proof} This is \cite[Ch.III, Prop.6.2, 8.2, 8.6]{Har}.
We check only (i). There is a natural commutative diagram of
ringed spaces :
$$
\xymatrix{ (X,\cO_{\!X}) \ar[r]^-f \ar[rd]_{\bar f} &
(Y,\cO_Y) \ar[r]^-g \ar[rd]_{\bar g} & (Z,\cO_{\!Z}) \\
& (Y,f_*\cO_{\!X}) \ar[rd]_\phi \ar[u]_\alpha & (Z,g_*\cO_Y) \ar[u]_\gamma \\
& & (Z,(g\circ f)_*\cO_{\!X}) \ar[u]_\beta
}$$
where :
\set\begin{equation}\label{eq_missing-arrow}
\phi\circ\bar f=\bar{g\circ f}.
\end{equation}
By claim \ref{cl_right-to-forget}, the (forgetful)
functors $\alpha_*$, $\beta_*$ and $\gamma_*$ admit
right adjoints :
$$
\begin{aligned}
\alpha^\flat & :\cO_Y\Mod\to f_*\cO_{\!X}\Mod & \quad : & \quad
\cF\mapsto\cHom_{\cO_Y}(f_*\cO_{\!X},\cF) \\
\beta^\flat & :g_*\cO_Y\Mod\to (g\circ f)_*\cO_{\!X}\Mod & \quad : & \quad
\cF\mapsto\cHom_{g_*\cO_Y}((g\circ f)_*\cO_{\!X},\cF) \\
\gamma^\flat &:\cO_{\!Z}\Mod\to g_*\cO_Y\Mod & \quad : & \quad
\cF\mapsto\cHom_{\cO_{\!Z}}(g_*\cO_Y,\cF).
\end{aligned}
$$
Likewise, $(\gamma\circ\beta)_*$ admits a right adjoint
$(\gamma\circ\beta)^\flat$ and the natural identification
$$
\gamma_*\circ\beta_*\isom(\gamma\circ\beta)_*
$$
induces a natural isomorphisms of functors :
\set\begin{equation}\label{eq_combo}
(\gamma\circ\beta)^\flat:=\cHom_{\cO_{\!Z}}((g\circ f)_*\cO_{\!X},-)
\isom\beta^\flat\circ\gamma^\flat.
\end{equation}
(Notice that, in general, there might be several choices of
such natural transformations \eqref{eq_combo}, but for
the proof of (iii) -- as well as for the construction
of $f^!$ for more general morphisms $f$ -- it is essential
to make an explicit and canonical choice.)

\begin{claim}\label{cl_compose}
(i)\ \ $\alpha_*\circ\phi^*=\bar g{}^*\circ\beta_*$.
\begin{enumerate}
\addenu
\item
The natural commutative diagram of sheaves :
$$
\xymatrix{ \bar g{}^{-1}g_*\cO_Y \ar[r] \ar[d] & \cO_Y \ar[d] \\
\bar g{}^{-1}(g\circ f)_*\cO_{\!X} \ar[r] & f_*\cO_{\!X}
}$$
is cocartesian.
\end{enumerate}
\end{claim}
\begin{pfclaim} (i) is an easy consequence of (ii). Assertion (ii)
can be checked on the stalks, hence we may assume that $Z=\Spec\,A$,
$Y=\Spec\,B$ and $X=\Spec\,C$ are affine schemes. Let
$\fp\in Y$ be any prime ideal, and set $\fq:=g(\fp)\in Z$;
then $(\bar g{}^{-1}g_*\cO_Y)_\fp=B_\fq$ (the $A_\fq$-module
obtained by localizing the $A$-module $B$ at the prime $\fq$).
Likewise, $(\bar g{}^{-1}(g\circ f)_*\cO_{\!X})_\fp=C_\fq$ and
$(f_*\cO_{\!X})_\fp=C_\fp$. Hence the claim boils down to the
standard isomorphism : $C_\fq\otimes_{B_\fq}B_\fp\simeq C_\fp$.
\end{pfclaim}

Let $I^\bullet$ be a bounded below
complex of injective $\cO_{\!Z}$-modules; by lemma \ref{lem_triv-dual}(i),
$\gamma^\flat I^\bullet$ is a complex of injective $g_*\cO_Y$-modules;
taking into account \eqref{eq_combo} we deduce a natural isomorphism
of functors on $\sD^+(\cO_{\!Z}\Mod)$ :
\set\begin{equation}\label{eq_der-combo}
R(\gamma\circ\beta)^\flat\isom R(\beta^\flat\circ\gamma^\flat)=
R\beta^\flat\circ R\gamma^\flat.
\end{equation}
Furthermore, \eqref{eq_missing-arrow} implies that
$\bar{g\circ f}{}^*=\bar f{}^*\circ\phi^*$. Combining with
\eqref{eq_der-combo}, we see that the sought $\xi_{f,g}$ is
a natural transformation :
$$
\bar f{}^*\circ\phi^*\circ R\beta^\flat\circ R\gamma^\flat\to
\bar f{}^*\circ R\alpha^\flat\circ\bar g{}^*\circ R\gamma^\flat
$$
of functors on $\sD^+(\cO_{\!Z}\Mod)$. Hence (i) will follow from :

\begin{claim}\label{cl_need-later}
(i)\ \ There exists a natural isomorphism of functors :
$$
g_*\cO_Y\Mod\to f_*\cO_{\!X}\Mod\quad :\quad
\phi^*\circ\beta^\flat\isom\alpha^\flat\circ\bar g{}^*.
$$
\begin{enumerate}
\addenu
\item
The natural transformation :
$$
R(\alpha^\flat\circ\bar g{}^*)\to
R\alpha^\flat\circ\bar g{}^*
$$
is an isomorphism of functors :
$\sD^+(g_*\cO_Y\Mod)\to\sD^+(f_*\cO_{\!X}\Mod)$.
\end{enumerate}
\end{claim}
\begin{pfclaim}[] First we show how to construct a natural
transformation as in (i); this is the same as exhibiting a
map of functors :
$\alpha_*\circ\phi^*\circ\beta^\flat\to\bar g{}^*$.
In view of claim \ref{cl_compose}(i), the latter can
be defined as the composition of $\bar g{}^*$ and
the counit of adjunction
$\beta_*\circ\beta^\flat\to\one_{g_*\cO_Y\Mod}$.

Next, recall the natural transformation :
\set\begin{equation}\label{eq_stupid}
\bar g{}^*\cHom_{g_*\cO_Y}(\cF,\cG)\to
\cHom_{\cO_Y}(\bar g{}^*\cF,\bar g{}^*\cG)
\quad\text{for every $g_*\cO_Y$-modules $\cF$ and $\cG$}
\end{equation}
(defined so as to induce the pull-back map on global $\Hom$ functors).
Notice that the natural map $\bar g{}^*\bar g_*\cF\to\cF$
(counit of adjunction) is an isomorphism for every quasi-coherent
$\cO_Y$-module $\cF$. Especially, if $\cF=g_*\cA$ for a
quasi-coherent $\cO_Y$-module $\cA$, then \eqref{eq_stupid}
takes the form :
\set\begin{equation}\label{eq_stupider}
\bar g{}^*\cHom_{g_*\cO_Y}(\bar g_*\cA,\cG)\to
\cHom_{\cO_Y}(\cA,\bar g{}^*\cG).
\end{equation}
Furthermore, by inspecting the definitions (and the proof of
claim \ref{cl_right-to-forget}), one verifies easily that the
map $\phi^*\circ\beta^\flat(\cG)\to\alpha^\flat\circ\bar g{}^*(\cG)$
constructed above is the same as the map \eqref{eq_stupider},
taken with $\cA=f_*\cO_{\!X}$.
Notice that, since $f$ is pseudo-coherent, $f_*\cO_{\!X}$ is even
a finitely presented $\cO_Y$-module; thus, in order to conclude
the proof of (i), it suffices to show that \eqref{eq_stupider}
is an isomorphism for every finitely presented $\cO_Y$-module $\cA$
and every $g_*\cO_Y$-module $\cG$. To this aim, we may assume
that $Z$ and $Y$ are affine; then we may find a presentation
$\underline\cE:=(\cO_Y^{\oplus p}\to\cO_Y^{\oplus q}\to\cA\to 0)$.
We apply the natural transformation \eqref{eq_stupider} to
$\underline\cE$, thereby obtaining a commutative ladder with
left exact rows; then the five-lemma reduces the assertion to
the case where $\cA=\cO_Y$, which is obvious. To prove assertion
(ii) we may again assume that $Y$ and $Z$ are affine, in which
case we can find a resolution $L_\bullet\to f_*\cO_{\!X}$ consisting
of free $\cO_Y$-modules of finite rank. For a given bounded below
complex $I^\bullet$ of injective $g_*\cO_Y$-modules, choose a
resolution $\bar g{}^*I^\bullet\to J^\bullet$ consisting of
injective $\cO_Y$-modules. We deduce a commutative ladder :
$$
\xymatrix{ \bar g{}^*\cHom^\bullet_{g_*\cO_Y}(\bar g_*\circ f_*\cO_{\!X},I^\bullet)
\ar[r]^-{\mu_1} \ar[d]_{\lambda_1} &
\cHom^\bullet_{\cO_Y}(f_*\cO_{\!X},\bar g{}^*I^\bullet) \ar[r]^-{\mu_3} \ar[d] &
\cHom^\bullet_{\cO_Y}(f_*\cO_{\!X},J^\bullet) \ar[d]^{\lambda_2} \\
\bar g{}^*\cHom^\bullet_{g_*\cO_Y}(\bar g_*L_\bullet,I^\bullet) \ar[r]^-{\mu_2} &
\cHom^\bullet_{\cO_Y}(L_\bullet,\bar g{}^*I^\bullet) \ar[r]^-{\mu_4} &
\cHom^\bullet_{\cO_Y}(L_\bullet,J^\bullet).
}$$
Since $\bar g$ is flat, $\lambda_1$ is a quasi-isomorphism, and the same
holds for $\lambda_2$. The maps $\mu_1$ and $\mu_2$ are of the form
\eqref{eq_stupider}, hence they are quasi-isomorphisms, by the foregoing
proof of (i). Finally, it is clear that $\mu_4$ is a quasi-isomorphism
as well, hence the same holds for $\mu_3\circ\mu_1$, and (ii) follows.
\end{pfclaim}
\end{proof}

\sset\subsubsection{}\label{subsec_sharp-flat}
Let now $f:X\to Y$ be an {\em embeddable} morphism,
{\em i.e.} such that it can be factored as a composition
$f=g\circ h$ where $h:X\to Z$ is a finite pseudo-coherent morphism,
$g:Z\to Y$ is smooth and separated, and the fibres of
$g$ have bounded dimension. One defines :
$$
f^!:=h^\flat\circ g^\sharp.
$$

\begin{lemma}\label{lem_compat-for_!} Let
$X\xrightarrow{f}Y\xrightarrow{g}Z\xrightarrow{h}W$
be three embeddable morphisms of schemes, such that the
compositions $g\circ f$, $h\circ g$ and $h\circ g\circ f$
are also embeddable. Then we have :
\begin{enumerate}
\item
There exists a natural isomorphism of functors
$$
\psi_{g,f}:(g\circ f)^!\isom f^!\circ g^!.
$$
Especially, $f^!$ is independent -- up to a
natural isomorphism of functors -- of the choice of factorization
as a finite morphism followed by a smooth morphism.
\item
The diagram of functors :
$$
\xymatrix{ (h\circ g\circ f)^! \ar[rr]^{\psi_{h,g\circ f}}
\ar[d]_{\psi_{h\circ g,f}} & &
(g\circ f)^!\circ h^! \ar[d]^{\psi_{g,f}\circ h^!} \\
f^!\circ(h\circ g)^! \ar[rr]^{f^!(\psi_{h,g})} & & f^!\circ g^!\circ h^!
}$$
commutes.
\end{enumerate}
\end{lemma}
\begin{proof} The proof is a complicated verification, starting
from proposition \ref{prop_sharp-flat}. See
\cite[Ch.III, Th.8.7]{Har} for details.
\end{proof}

\sset\subsubsection{}\label{subsec_wants-to-commute} Let
$X\xrightarrow{f}Y\xrightarrow{g}Z$ be two finite and
pseudo-coherent morphisms of schemes. The map $f^\natural:\cO_Y\to
f_*\cO_{\!X}$ induces a natural transformation :
\set\begin{equation}\label{eq_exhibit}
\cHom^\bullet_{\cO_{\!Z}}(g_*(f^\natural),-):
\beta_*\circ(\gamma\circ\beta)^\flat\to\gamma^\flat
\end{equation}
(notation of the proof of proposition \ref{prop_sharp-flat}(i))
and we wish to conclude this section by exhibiting another
compatibility involving \eqref{eq_exhibit} and the isomorphism
$\xi_{f,g}$. To this aim, we consider the composition of isomorphism
of functors on $\sD^+(\cO_{\!Z}\Mod)_\qcoh$ :
\set\begin{equation}\label{eq_boxes}
Rf_*\circ\bar{g\circ f}{}^*
\xrightarrow{\xymatrix{*+[o][F]{\text{\tiny 1}}}}
\alpha_*\circ R\bar f_{\!*}\circ\bar f{}^*\circ\phi^*
\xrightarrow{\xymatrix{*+[o][F]{\text{\tiny 2}}}}
\alpha_*\circ\phi^*
\xrightarrow{\xymatrix{*+[o][F]{\text{\tiny 3}}}}
\bar g{}^*\circ\beta_*
\end{equation}
where :

\begin{itemize}
\item[$\xymatrix{*+[o][F]{1}}$] is deduced from \eqref{eq_missing-arrow} and
the decomposition $f=\alpha\circ\bar f$.
\item[$\xymatrix{*+[o][F]{2}}$] is induced from the counit of adjunction
$\bar\eps:R\bar f_{\!*}\circ\bar f{}^*\to\one$ provided by lemma
\ref{lem_Groth-dual}.
\item[$\xymatrix{*+[o][F]{3}}$] is the identification of claim
\ref{cl_compose}.
\end{itemize}
We define a morphism of functors as a composition :
\set\begin{equation}\label{eq_new-boxen}
Rf_*\circ(g\circ f)^\flat=
Rf_*\circ\bar{g\circ f}{}^*\circ R(\gamma\circ\beta)^\flat
\xrightarrow{\xymatrix{*+[o][F]{\scriptstyle 4}}}
\bar g{}^*\circ\beta_*\circ R(\gamma\circ\beta)^\flat
\xrightarrow{\xymatrix{*+[o][F]{\scriptstyle 5}}}
g^\flat
\end{equation}
where :
\begin{itemize}
\item[$\xymatrix{*+[o][F]{4}}$] is the isomorphism
$\eqref{eq_boxes}\circ R(\gamma\circ\beta)^\flat$.
\item[$\xymatrix{*+[o][F]{5}}$] is the composition of
$\bar g{}^*$ and the morphism
$\beta_*\circ R(\gamma\circ\beta)^\flat\to R\gamma^\flat$
deduced from $\eqref{eq_exhibit}$.
\end{itemize}
We have a natural diagram of functors on $\sD^+(\cO_{\!Z}\Mod)_\qcoh$ :
\set\begin{equation}\label{eq_wants-to-commute}
{\diagram R(g\circ f)_*\circ(g\circ f)^\flat \ar[r]^-*+[o][F]{\text{\tiny 6}}
\ar[d]_{R(g\circ f)_*(\xi_{f,g})\ } & Rg_*\circ Rf_*\circ(g\circ f)^\flat
\ar[r]^-*+[o][F]{\text{\tiny 9}}  \ar[d]^{Rg_*\circ Rf_*(\xi_{f,g})\ } &
Rg_*\circ g^\flat \\
R(g\circ f)_*\circ(f^\flat\circ g^\flat) \ar[r]^-*+[o][F]{\text{\tiny 7}} &
Rg_*\circ Rf_*\circ(f^\flat\circ g^\flat) \rdouble &
Rg_*\circ(Rf_*\circ f^\flat)\circ g^\flat \ar[u]_*+[o][F]{\text{\tiny 8}}
\enddiagram}\end{equation}
where :
\begin{itemize}
\item[$\xymatrix{*+[o][F]{6}}$] and $\xymatrix{*+[o][F]{7}}$ are induced
by the natural isomorphism $R(g\circ f)_*\isom Rg_*\circ Rf_*$.
\item[$\xymatrix{*+[o][F]{8}}$] is induced by the counit of the adjunction
$\eps:Rf_*\circ f^\flat\to\one$.
\item[$\xymatrix{*+[o][F]{9}}$] is $Rg_*\circ\eqref{eq_new-boxen}$.
\end{itemize}

\begin{lemma}\label{lem_wants-to-commute}
In the situation of \eqref{subsec_wants-to-commute},
diagram \eqref{eq_wants-to-commute} commutes.
\end{lemma}
\begin{proof} The diagram splits into left and right subdiagrams,
and clearly the left subdiagram commutes, hence it remains to
show that the right one does too; to this aim, it suffices to consider
the simpler diagram :
$$
\xymatrix{
Rf_*\circ(g\circ f)^\flat \ar[r] \ar[d]_{Rf_*(\xi_{f,g})} & g^\flat \\
Rf_*\circ f^\flat\circ g^\flat \ar[ur]
}$$
whose horizontal arrow is \eqref{eq_new-boxen}, and whose upward arrow
is deduced from the counit $\eps$. However, the counit $\eps$, used in
$\xymatrix{*+[o][F]{8}}$, can be expressed in terms of the counit
$\bar\eps:R\bar f_*\circ\bar f{}^*\to\one$, used in $\xymatrix{*+[o][F]{2}}$,
therefore we are reduced to considering the diagram of functors on
$\sD^+(\cO_{\!Z}\Mod)_\qcoh$ :
$$
\xymatrix{
\alpha_*\circ R\bar f_*\circ\bar f{}^*\circ\phi^*\circ
R(\gamma\circ\beta)^\flat \ar[r]^-*+[o][F]{\scriptstyle a}
\ar[dd]_{Rf_*(\xi_{f,g})} &
\alpha_*\circ\phi^*\circ R(\gamma\circ\beta)^\flat
\ar[d]_*+[o][F]{\scriptstyle c} \ar[r]^-*=<11pt>[o][F]{\scriptstyle e} &
\bar g{}^*\circ\beta_*\circ R(\gamma\circ\beta)^\flat
\ar[dd]^*=<11pt>[o][F]{\scriptstyle 5} \\
& \alpha_*\circ\phi^*\circ R\beta^\flat\circ R\gamma^\flat
\ar[d]_-*+[o][F]{\scriptstyle d} \\
\alpha_*\circ R\bar f_*\circ\bar f{}^*\circ R\alpha^\flat\circ g^\flat
\ar[r]^-*+[o][F]{\scriptstyle b} &
\alpha_*\circ R\alpha^\flat\circ\bar g{}^*\circ R\gamma^\flat
\ar[r]^-*+[o][F]{\scriptstyle f} &
g^\flat=\bar g{}^*\circ R\gamma^\flat \\
}$$
where :
\begin{itemize}
\item[$\xymatrix{*=<14pt>[o][F]{a}}$] and $\xymatrix{*=<14pt>[o][F]{b}}$
are induced from $\bar\eps$.
\item[$\xymatrix{*=<14pt>[o][F]{c}}$] is induced by the
isomorphism \eqref{eq_der-combo}.
\item[$\xymatrix{*=<14pt>[o][F]{d}}$] is induced
by the isomorphism of claim \ref{cl_need-later}(i).
\item[$\xymatrix{*=<14pt>[o][F]{e}}$]
is induced by the identity of claim \ref{cl_compose}.
\item[$\xymatrix{*=<14pt>[o][F]{f}}$] is induced by the counit
$\alpha_*\circ R\alpha^\flat\to\one$.
\end{itemize}
Now, by inspecting the construction of $\xi_{f,g}$ one checks that the
left subdiagram of the latter diagram commutes; moreover, from claim
\ref{cl_need-later}(ii) (and from lemma \ref{lem_triv-dual}(i)), one
sees that the right subdiagram is obtained by evaluating on complexes
of injective modules the analogous diagram of functors on $\cO_{\!Z}\Mod$ :
$$
\xymatrix{
\alpha_*\circ\phi^*\circ(\gamma\circ\beta)^\flat
\ar[r]^-*=<11pt>[F]{\scriptstyle e} \ar[d]_-*=<11pt>[F]{\scriptstyle c} &
\bar g{}^*\circ\beta_*\circ(\gamma\circ\beta)^\flat
\ar[dd]^*=<11pt>[F]{\scriptstyle 5} \\
\alpha_*\circ\phi^*\circ\beta^\flat\circ\gamma^\flat
\ar[d]_-*=<11pt>[F]{\scriptstyle d} \\
\alpha_*\circ\alpha^\flat\circ\bar g{}^*\circ\gamma^\flat
\ar[r]^-*=<11pt>[F]{\scriptstyle f} &
\bar g^*\circ\gamma^\flat.
}$$
Furthermore, by inspecting the proof of claim \ref{cl_need-later}(i),
one verifies that $\xymatrix{*=<14pt>[F]{f}}\circ\xymatrix{*=<14pt>[F]{d}}$
is the same as the composition
$$
\alpha_*\circ\phi^*\circ\beta^\flat\circ\gamma^\flat
\xrightarrow{\xymatrix{*+[o][F]{\scriptstyle i}}}
\bar g{}^*\circ\beta_*\circ\beta^\flat\circ\gamma^\flat
\xrightarrow{\xymatrix{*+[o][F]{\scriptstyle j}}}
\bar g{}^*\circ\gamma^\flat
$$
where $\xymatrix{*=<14pt>[o][F]{j}}$ is deduced from the counit
$\beta_*\circ\beta^\flat\to\one$, and $\xymatrix{*=<14pt>[o][F]{i}}$
is deduced from the identity of claim \ref{cl_compose}.
Hence we come down to showing that the diagram :
$$
\xymatrix{
\beta_*\circ(\gamma\circ\beta)^\flat\ar[d]_{\beta_*\eqref{eq_combo}}
\ar[rr]^-{\eqref{eq_exhibit}} & & \gamma^\flat \\
\beta_*\circ\beta^\flat\circ\gamma^\flat
\ar[rru]_-*+[o][F]{\scriptstyle k}
}$$
commutes in $g_*\cO_Y\Mod$, where $\xymatrix{*=<14pt>[o][F]{k}}$
is deduced from the counit $\beta_*\circ\beta^\flat\to\one$. This is
an easy verification, that shall be left to the reader.
\end{proof}

\subsection{Cousin complexes}\label{sec_Cousin}
In this section we study a special type of complexes of
abelian sheaves, called {\em Cousin complexes}, defined
on a large class of topological spaces. These
complexes play an important role in the further development
of duality theory for locally coherent schemes : in many
cases the techniques explained in this section will
allow to exhibit canonical representatives for the dualizing
complexes that shall be investigated at length in section
\ref{subsec_duals}.

\sset\subsubsection{}\label{subsec_here}
We resume the notation of \eqref{subsec_loc-closed-supp},
so $(X,\cA)$ is a ringed topological space, $\phi:V\to X$ a
locally closed immersion, we denote by $\bar V$ the topological
closure of $V$ in $X$, and we set $\partial V:=\bar V\!\setminus\!V$.
In \eqref{subsec_Gamma_Z} we have introduced the functor
$\Gamma_{\!\bar V}$ of sections with support in $\bar V$,
and its sheaf-valued variant $\underline\Gamma_{\bar V}$.
We now extend these definitions to the case of a general
locally closed immersion, by setting
$$
\underline\Gamma_V:=
\phi_*\circ\phi^*\circ\underline\Gamma_{\bar V}:
\cA\Mod\to\cA\Mod
\qquad
\Gamma_V:=\Gamma\circ\underline\Gamma_V:\cA\Mod\to\cA(X)\Mod.
$$
The functor $\underline\Gamma_V$ admits a derived functor
$$
R\underline\Gamma_V:\sD^+(\cA\Mod)\to\sD^+(\cA\Mod)
$$
obtained as usual by evaluating on injective resolutions
of given bounded below complexes.

\begin{remark}\label{rem_consistence}
(i)\ \
We have the following explicit description of 
$\underline\Gamma_V$. Let $k:V\to X\!\setminus\!\partial V$
and $l:X\!\setminus\!\partial V\to X$ be respectively the
closed and open immersions. Since $\phi=l\circ k$, there
exist natural isomorphisms
$$
\underline\Gamma_V\isom
\phi_*\circ k^*\circ l^*\circ\underline\Gamma_{\bar V}\isom
l_*\circ k_*\circ k^*\circ\underline\Gamma_V\circ l^*\isom
l_*\circ\underline\Gamma_V\circ l^*
$$
where the last map is induced by the natural identification
$k_*\circ k^*\circ\underline\Gamma_V\isom\underline\Gamma_V$
provided by lemma \ref{lem_upside-down}. There follows
a natural isomorphism
$$
\underline\Gamma_V\cF(S)\isom
\underline\Gamma_{\bar V}\cF(S\!\setminus\!\partial V)
\qquad
\text{for every open subset $S$ of $X$}.
$$

(ii)\ \
Notice that $\partial V=\emptyset$ if $\phi$ is a closed
immersion; in this case, (i) says that the functor defined
in \eqref{subsec_here} is naturally isomorphic
to its namesake introduced in \eqref{subsec_Gamma_Z}, so
the current notation is consistent with that of section
\ref{subsec_supp}.

(iii)\ \
For a general locally closed subset $V$ and any
$\cA$-module $\cF$, the natural monomorphism
$\iota_\cF:\underline\Gamma_{\bar V}\cF\to\cF$ induces
a natural transformation
\set\begin{equation}\label{eq_always-mono}
\underline\Gamma_V\cF\to\phi_*\phi^*\cF
\end{equation}
which is a monomorphism for every such $\cF$, since
the functor $\phi_*\phi^*$ is left exact. Notice that
if $V$ is open in $X$, the map $\phi^*\iota_\cF$ is an
isomorphism, so the same holds for \eqref{eq_always-mono}.
\end{remark}

\begin{proposition}\label{prop_canonic-adjun}
In the situation of \eqref{subsec_here},
the functor $\underline\Gamma_V$ is right adjoint to $t_V$.
\end{proposition}
\begin{proof} For $\phi$ a closed immersion, remark
\ref{rem_consistence}(ii) reduces the assertion to lemma
\ref{lem_adj-Gamma_Z}(i), since in this case $\phi_!=\phi_*$.
The assertion is also clear in case $\phi$ is an open
immersion, since in this case $\underline\Gamma_V$ is
naturally isomorphic to the functor $\phi_*\circ\phi^*$,
by remark \ref{rem_consistence}(iii). Notice that
$V=\bar V\cap(X\!\setminus\!\partial V)$. By proposition
\ref{prop_double-truncate}(i) (and by remark
\ref{rem_adjoint-transf}(i)) we are reduced
to exhibiting a natural isomorphism of functors
$$
\underline\Gamma_V\isom
\underline\Gamma_{X\setminus\partial V}\circ
\underline\Gamma_{\bar V}.
$$
But the existence of such an isomorphism is asserted
already in remark \ref{rem_consistence}(i).
\end{proof}

\sset\subsubsection{}\label{subsec_can-discussion}
For our discussion, it will be important to know that
in fact there exists a canonical adjunction for the
pair $(t_V,\underline\Gamma_V)$. In view of the proof
of proposition \ref{prop_canonic-adjun} (and of
remark \ref{rem_adjoint-transf}(i)), it suffices
to exhibit a canonical pair $(\eta^V,\eps^V)$ of units
and counits for this adjunction, in case $\phi$ is
either an open or a closed immersion.

$\bullet$\ \
Suppose first that $\phi$ is a closed immersion. Then
we want natural transformations
$$
\eta^V_\cF:\cF\to\underline\Gamma_V\circ\phi_*\circ\phi^*\cF
\qquad
\eps^V_\cF:\phi_*\circ\phi^*\circ\underline\Gamma_V\cF\to\cF
$$
related by the triangular identities of \eqref{subsec_adj-pair}.
However, remark \ref{rem_consistence}(ii) yields a
natural isomorphism
$\eta_1:\underline\Gamma_V\phi_*\phi^*\cF\isom\phi_*\phi^*\cF$,
and on the other hand, we have a canonical unit
$\eta_2:\cF\to\phi_*\phi^*\cF$ for the adjunction
$(\phi^*,\phi_*)$, so we take
$\eta^V_\cF:=\eta_1^{-1}\circ\eta_2$. Likewise, let
$\eps_1:\underline\Gamma_V\cF\to\cF$ be the natural
monomorphism; we take $\eps^V_\cF:=\eps_1\circ f_\cF$,
where
$f_\cF:t_V\circ\underline\Gamma_V\cF\isom\underline\Gamma_V\cF$
is the isomorphism provided by lemma \ref{lem_upside-down}.
The triangular identities may then be checked on the
stalks : we leave the details to the reader.

$\bullet$\ \
Next, if $\phi$ is an open immersion, we need to find
natural transformations
$$
\eta^V_\cF:\cF\to\phi_*\circ\phi^*\circ\phi_!\circ\phi^*\cF
\qquad
\eps^V_\cF:\phi_!\circ\phi^*\circ\phi_*\circ\phi^*\cF\to\cF
$$
again related by the same triangular identities. We take
for $\eta^V_\cF$ (resp. $\eps^V_\cF$) the composition of
the natural unit $\cF\to\phi_*\phi^*\cF$ of the adjunction
$(\phi^*,\phi^*)$ (resp. the counit $\phi_!\circ\phi^*\cF\to\cF$
of the adjunction $(\phi_!,\phi^*)$), with the isomorphism
$$
\phi_*\phi^*\cF\isom\phi_*\phi^*\phi_!\phi^*\cF
\qquad
\text{(\ resp.
$\phi_!\phi^*\phi_*\phi^*\cF\isom\phi_!\phi^*\cF$\ )}
$$
deduced from the natural identification of $\phi^*\circ\phi_!$
(resp. of $\phi^*\circ\phi_*$) with the identity functor
of $\phi^*\cA\Mod$. The triangular identities are as usual
checked on the stalks.

\sset\subsubsection{}\label{subsec_variant-i-j}
Due to proposition \ref{prop_canonic-adjun}, and the
subsequent observations in \eqref{subsec_can-discussion}, all
the constructions of proposition \ref{prop_double-truncate}
have natural counterparts for the functors $\underline\Gamma_V$.
Namely :

$\bullet$\ \
If $V$ and $V'$ are any two locally closed subsets of $X$,
there is a natural isomorphism
$$
\beta_{V,V'}:
\underline\Gamma_V\circ\underline\Gamma_{V'}\isom
\underline\Gamma_{V\cap V'}
$$
and if $V''\subset X$ is any other locally closed subset,
the diagram of functors
\set\begin{equation}\label{eq_triple-compos}
{\diagram
\underline\Gamma_V\circ\underline\Gamma_{V'}\circ
\underline\Gamma_{V''}
\ar[rrr]^-{\underline\Gamma_V*\beta_{V',V''}}
\ar[d]_{\beta_{V,V'}*\underline\Gamma_{V''}} & & &
\underline\Gamma_V\circ\underline\Gamma_{V'\cap V''}
\ar[d]^{\beta_{V,V'\cap V''}}  \\
\underline\Gamma_{V\cap V'}\circ\underline\Gamma_{V''}
\ar[rrr]^-{\beta_{V\cap V',V''}} & & &
\underline\Gamma_{V\cap V'\cap V''}
\enddiagram}
\end{equation}
commutes.

$\bullet$\ \
Moreover, suppose that the pair $(V,V')$ satisfies
either of the conditions (a), (b) of proposition
\ref{prop_double-truncate}(iii); then there exists
a natural morphism
$$
\gamma_{V,V'}:\underline\Gamma_V\to\underline\Gamma_{V'}.
$$
which corresponds to $c^{V,V'}$ under the natural
adjunction described in \eqref{subsec_can-discussion},
and if $V''$ is any third locally closed subset
of $X$, we get a commutative diagram
$$
\xymatrix{
\underline\Gamma_{V''}\circ\underline\Gamma_V
\ar[rrr]^-{\underline\Gamma_{V''}*\gamma_{V,V'}}
\ar[d]_{\beta_{V'',V}} & & &
\underline\Gamma_{V''}\circ\underline\Gamma_{V'}
\ar[d]^{\beta_{V'',V'}} \\
\underline\Gamma_{V''\cap V}
\ar[rrr]^-{\gamma_{V''\cap V,V''\cap V'}} & & &
\underline\Gamma_{V''\cap V'}.
}$$

$\bullet$\ \
Lastly, in the situation of proposition
\ref{prop_double-truncate}(v) we have the identity :
$$
\gamma_{V',V''}\circ\gamma_{V,V'}=\gamma_{V,V''}:
\underline\Gamma_V\to\underline\Gamma_{V''}.
$$

\begin{remark}\label{rem_different-def}
Let $X$ be a topological space, $j:U\to X$ an open
immersion, $\phi:V\to X$ and $\phi':V'\to X$ two
locally closed immersions fulfilling either of conditions
(a),(b) of proposition \ref{prop_double-truncate}(iii).
The transformations $\gamma_{V,V'}$ can also be described
explicitly, thanks to the following considerations.

(i)\ \
First, notice that the functors $t_\bullet$ on
$\cA$-modules, the unit $\eta^\bullet$ and the counit
$\eps^\bullet$ are all compatible with restriction to
$U$, {\em i.e.} we have identities of functors on
$j^*\cA$-modules and of morphisms between such functors :
$$
j^**t_V=t_{U\cap V}
\quad
j^**\eta^V=\eps^{V\cap U}
\quad
j^**\eps^V=\eta^{V\cap U}
\qquad
\text{for every locally closed subset $V\subset X$}.
$$
This is obvious for the truncation functors, and for
the units and counits it follows by a simple inspection
of the definition. Then, the uniqueness properties of
the transformations $\tau^{\bullet\bullet}$ and
$c^{\bullet\bullet}$ imply that the same compatibility
with restriction to $U$ holds for the latter.

(ii)\ \
It follows formally from (i) that the isomorphisms
$\beta_{\bullet\bullet}$ and the transformations
$\gamma_{\bullet\bullet}$ are likewise compatible
with restriction to $U$. On the other hand, notice
that the restriction map
$$
\Hom_{\cA\Mod}(t_{V'}\cF,\cG)\to
\Hom_{j^*\cA\Mod}((t_{V'}\cF)_{|U},\cG_{|U})
$$
is a bijection for every $\cF,\cG\in\Ob(\cA\Mod)$,
whenever $V'\subset U$. We conclude that the restriction map
$$
\Hom_{\cA\Mod}(\cF,\underline\Gamma_{V'}\cG)\to
\Hom_{j^*\cA\Mod}(\cF_{|U},\underline\Gamma_{V'\cap U}\cG_{|U})
$$
is also bijective when $V'\subset U$. Letting
$\cF:=\underline\Gamma_V\cG$, we deduce that
$\gamma_{V,V'}$ is determined by its restriction
to $U$, which is just $\gamma_{V\cap U,V'\cap U}$.

(iii)\ \
Now, suppose first that $\psi:V\to V'$ is a closed
immersion, and denote by $\bar V$ and $\bar V{}'$ the
topological closures in $X$ of $V$ and respectively
$V'$; notice that the unit of the adjunction
$(\psi^*,\psi_*)$ induces an isomorphism
\set\begin{equation}\label{eq_towards-ch-supps}
\phi'{}^*\circ\Gamma_{\bar V}\isom
\psi_*\circ\psi^*\circ\phi'{}^*\circ\Gamma_{\bar V}
\end{equation}
since $(\phi'{}^*\circ\Gamma_{\bar V}\cF)_{|V'\setminus V}=0$
for every $\cA$-module $\cF$ (lemma \ref{lem_upside-down}).
We claim that $\gamma_{V,V'}$ is the composition
$$
\phi_*\circ\phi^*\circ\underline\Gamma_{\bar V}\isom
\phi'_*\circ\psi_*\circ\psi^*\circ\phi'{}^*\circ
\underline\Gamma_{\bar V}\isom
\phi'_*\circ\phi'{}^*\circ\underline\Gamma_{\bar V}\to
\phi'_*\circ\phi'{}^*\circ\underline\Gamma_{\bar V{}'}.
$$
where the second map is the inverse of
\eqref{eq_towards-ch-supps} and the third is deduced
from $\gamma_{\bar V,\bar V'}:
\underline\Gamma_{\bar V}\to\underline\Gamma_{\bar V{}'}$.
Indeed, notice that the foregoing construction is
compatible with restriction to $U$; taking into account
(iii), we may then replace $X$ by its open subset
$X\setminus\partial V'$, and assume that both $\phi$
and $\phi'$ are closed immersions. In this case, the
assertion follows by direct inspection of the definition
of the unit and counit of the adjunction
$(t_\bullet,\underline\Gamma_\bullet)$, and of the
definition of $t_\bullet$, in the case of closed immersions.

(iv)\ \
Likewise, if $V'$ is an open subset of $V$, pick a
factorization of $\phi:V\to X$ as the composition of
a closed immersion $i:V\to U$ and an open immersion
$j:U\to X$. There is a unique open subset $U'\subset U$
such that $V'=U'\cap V$; denote by $j':U'\to X$ the
open immersion, and for any given $\cA$-module $\cF$,
let $\eta_\cF:\cF\to j'_*j'{}^*\cF$ be the unit of
adjunction; according to remark \ref{rem_consistence},
we may regard $\eta_\cF$ as a natural transformation
$$
\eta_\cF:\cF\to\underline\Gamma_{U'}\cF
$$
and we claim that $\gamma_{V,V'}$ is the composition
$$
\underline\Gamma_V
\xrightarrow{\ \underline\Gamma_V*\eta\ }
\underline\Gamma_V\circ\underline\Gamma_{U'}
\xrightarrow{\ \beta_{V,U'}\ }\underline\Gamma_{V'}.
$$
For the proof, we may argue as in (iii) to see that
$\gamma_{V,V'}$ is determined by its restriction to
$U'$, and since the proposed construction is compatible
with restriction to $U'$, we are then reduced to the
case where $U'=U=X$ and therefore $V=V'$. However,
clearly $c^{V,V}$ is the identity endomorphism of
$t_V$, so it remains only to prove that the proposed
map is the identity of $\underline\Gamma_V$, if
$V=V'$. The latter can be checked easily, whence
the claim. We could have used the above construction
as a definition for $\gamma_{V,V'}$ in the case of
open immersions, but the method we followed has the
advantage of dispensing us from verifying that the
resulting transformation is independent of the chosen
factorization of $\phi$.
\end{remark}

\sset\subsubsection{}\label{subsec_formalities}
Let $V\subset V'$ be a closed immersion of locally
closed subsets of a topological space $X$; from
remark \ref{rem_different-def}(iii), it is clear
that $\gamma_{V,V'}$ is a monomorphism. Moreover,
we have :

\begin{lemma}\label{lem_supp-left-ex}
In the situation of \eqref{subsec_formalities}, the
sequence
$$
\Sigma_{V,V'}
\quad :\quad
0\to\underline\Gamma_V\cF\xrightarrow{\ \gamma_{V,V'}\ }
\underline\Gamma_{V'}\cF\xrightarrow{\ \gamma_{V',V'\setminus V}\ }
\underline\Gamma_{V'\setminus V}\cF
$$
is left exact for every $\cA$-module $\cF$.
\end{lemma}
\begin{proof} Indeed, say that $V'=U\cap Z$ for
some closed subset $Z$ and open subset $U$; then
$V=\bar V\cap U$, with $\bar V$ the topological
closure of $V$ in $X$, and in light of
\eqref{subsec_variant-i-j} we see that there is
a natural isomorphism
$\Sigma_{V,V'}\isom\underline\Gamma_U(\Sigma_{\bar V,Z})$.
However, $\underline\Gamma_U$ is a left exact functor
(since it admits a left adjoint), so it suffices to
check that $\Sigma_{\bar V,Z}$ is left exact. Thus, we
may assume from start that $V$ and $V'$ are closed in
$X$. By the same token, we see that
$\Sigma_{V,V'}\simeq\underline\Gamma_{V'}(\Sigma_{V,X})$,
so we may further assume that $V'=X$, in which case the
sequence reduces to the initial segment of
\eqref{eq_first-two-terms}.
\end{proof}

\begin{proposition}\label{prop_derive-undergamma}
With the notation of \eqref{subsec_here}, the
following holds :
\begin{enumerate}
\item
If $\cF$ is a flabby $\cA$-module, the same holds for
$\underline\Gamma_V\cF$.
\item
Every flabby $\cA$-module is $\underline\Gamma_V$-acyclic.
\item
If $\cB\to\cA$ is any morphism of sheaves of rings on $X$,
the functor $R\underline\Gamma_V$ commutes with the forgetful
functor $\sD^+(\cA\Mod)\to\sD^+(\cB\Mod)$.
\item
If $\phi':V'\to X$ is any other locally closed immersion,
there exists a natural isomorphism of functors :
$$
R\underline\Gamma_{V\cap V'}\isom
R\underline\Gamma_V\circ R\underline\Gamma_{V'}
\qquad
\sD^+(\cA\Mod)\to\sD^+(\cA\Mod).
$$
\item
Suppose that $X$ is locally spectral and quasi-separated,
and $V$ is a constructible locally closed subset of $X$,
then we have :
\begin{enumerate}
\item
Every qc-flabby $\cA$-module is $\underline\Gamma_V$-acyclic.
\item
For every $i\in\N$, the functor $R^i\underline\Gamma_V$
commutes with filtered colimits of $\cA$-modules.
\item
For every $x\in V$, denote by $j_x:X(x)\to X$ the inclusion
map (notation of definition {\em\ref{def_special}(iii)}).
Then there exists a natural isomorphism
$$
j_x^*R\underline\Gamma_V\isom
R\underline\Gamma_{V(x)}
\qquad
\text{in $\sD^+(j_x^*\cA\Mod)$}.
$$
\end{enumerate}
\end{enumerate}
\end{proposition}
\begin{proof}(i): By lemmata \ref{lem_adj-Gamma_Z}(ii)
and \ref{lem_flabby-Gamma_Z}(i), we know already that
both assertions (i) and (ii) hold in case $\phi$ is a
closed immersion; for the general case, we may then
assume that $\phi$ is an open immersion. However, in
this case both $\phi_*$ and $\phi^*$ take flabby sheaves
to flabby sheaves (lemmata \ref{lem_flabby-is-local} and
\ref{lem_franziska}(ii)), and flabby sheaves are acyclic
for both of these functors (lemma \ref{lem_franziska}(ii)),
so the assertion follows from remark \ref{rem_consistence}(iii).

(iii) is clear from (ii), as every injective $\cA$-module
is a flabby $\cB$-module (lemma \ref{lem_franziska}(v)).

(iv) follows formally from (i), (ii) and
\eqref{subsec_variant-i-j} : details left to the reader.

(v.a): Since the assertion is local on $X$, we may
assume that $X$ is spectral, in which case $V$ is the
intersection of a closed constructible subset and an
open constructible subset.
By (iv) and lemma \ref{lem_franziska}(iii), it then
suffices to consider separately the cases where $V$
is either open or closed in $X$. By lemma
\ref{lem_flabby-Gamma_Z}(iii.b) we know already
the assertions in case $V$ is a closed constructible
subset of $X$, and the case where $V$ is open
follows easily from lemma \ref{lem_franziska}(i,iv).

(v.b): We reduce again to the case where $V$ is the
intersection of a closed constructible subset $Z$ and
an open constructible subset $U$. Since -- by (ii) --
the functor $\underline\Gamma_Z$ transforms injective
$\cA$-modules into $\underline\Gamma_U$-acyclic
$\cA$-modules, we have a convergent spectral sequence
$$
R^i\underline\Gamma_U\circ R^j\underline\Gamma_Z\cF
\Rightarrow R^{i+j}\underline\Gamma_V\cF.
$$
It then suffices to consider separately the case where
$V$ is either open or closed. These cases follow
respectively from lemma \ref{lem_flabby-Gamma_Z}(iii.b)
and claim \ref{cl_gosh}.

(v.c): Arguing as in the proof of (v.b), we reduce
to consider separately the cases where $V$ is either
open or closed. For $V$ closed, the assertion is
corollary \ref{cor_better-late}. For $V$ open, the
assertion follows from proposition
\ref{prop_dir-im-and-colim}(ii) and remark
\ref{rem_consistence}(iii), by arguing as in the
proof of corollary \ref{cor_better-late} : details
left to the reader.
\end{proof}

\begin{remark}
(i)\ \ 
Using proposition \ref{prop_derive-undergamma}, the
various compatibilities that have been found so far
beween the functors $\underline\Gamma_V$ extend
straightforwardly to their derived extensions.
Especially, diagram \eqref{eq_triple-compos} still
commutes, after replacing the functors
$\underline\Gamma_\bullet$ by their derived extensions
everywhere, and the natural transformations
$\beta_{\bullet\bullet}$ by their derived versions, provided
by proposition \ref{prop_derive-undergamma}(iii).

(ii)\ \
Likewise, if $V\subset V'$ (resp. $V'\subset V$) is any
closed (resp. open) immersion of locally closed subsets
of $X$, by deriving the morphisms $\gamma_{V,V'}$ of
\eqref{subsec_variant-i-j} we get a natural transformation
$$
R\underline\Gamma_V\to R\underline\Gamma_{V'}
\qquad
\text{(\ resp. $R\underline\Gamma_{V'}\to R\underline\Gamma_V$\ )}
\qquad
\text{in $\sD^+(\cA\Mod)$}
$$
and the compatibilities between the functors
$\gamma_{\bullet\bullet}$ and $\beta_{\bullet\bullet}$
explicited in \eqref{subsec_variant-i-j} hold,
{\em mutatis mutandis}, also for the functors
$R\underline\Gamma_\bullet$.
\end{remark}

\begin{lemma} In the situation of \eqref{subsec_formalities},
let $\cF$ be any flabby $\cA$-module. Then, for every open
subset $U\subset X$, the sequence
$$
0\to\underline\Gamma_V\cF(U)\to
\underline\Gamma_{V'}\cF(U)\to
\underline\Gamma_{V'\setminus V}\cF(U)\to 0
$$
induced by the sequence $\Sigma_{V,V'}$ of lemma
{\em\ref{lem_supp-left-ex}}, is short exact.
\end{lemma}
\begin{proof} We may assume that $U=X$, and it
suffices to show that the map
$g:\Gamma_{\!V'}\cF\to\Gamma_{\!V'\setminus V}\cF$
deduced from $\gamma_{V',V'\setminus V}$ is surjective,
whenever $\cF$ is flabby. To this aim, write
$V'\setminus V=U'\cap V'$ for some open subset
$U'\subset X$; then $g$ is naturally identified
with the map
$$
\Gamma_{\!X}\circ\underline\Gamma_{\!V'}\cF\to
\Gamma_{\!U'}\circ\underline\Gamma_{V'}\cF
$$
deduced from $\gamma_{X,U'}$. But $\underline\Gamma_{V'}\cF$
is still flabby (proposition \ref{prop_derive-undergamma}(i))
so we are reduced to the case where $V'=X$, in which case
$g$ is naturally identified with the restriction map
$\cF(X)\to\cF(X\!\setminus\!V)$, whence the contention.
\end{proof}

\sset\subsubsection{}\label{subsec_family-of-supps}
Let $(X,\cA)$ be a ringed topological space.
A {\em family of supports} on $X$ is a set $\Phi$
consisting of closed subsets of $X$, such that :
\begin{itemize}
\item
$\emptyset\in\Phi$.
\item
If $Z\in\Phi$ and $Z'\subset Z$ is any closed subset,
then $Z'\in\Phi$ as well.
\item
If $Z_1,Z_2\in\Phi$, then $Z_1\cup Z_2\in\Phi$ as well.
\end{itemize}
Now, let $\Phi$ and $\Phi'$ be any two families of
supports on $X$; we let $\Phi|\Phi'$ be the set whose
elements are all the pairs $\underline Z:=(Z,Z')$, where
$Z$ and $Z'$ are arbitrary elements of $\Phi$ and
respectively $\Phi'$, such that $Z'\subset Z$. We endow
$\Phi|\Phi'$ with a partial ordering, by declaring that,
for given pairs
$\underline Z_1:=(Z_1,Z'_1),\underline Z_2:=(Z_2,Z'_2)$
of $\Phi|\Phi'$, we have $\underline Z_1\leq\underline Z_2$
if $Z_1\subset Z_2$ and $Z'_1\subset Z'_2$. It is easily
seen that $\Phi|\Phi'$ is filtered.
To every object $\underline Z\in\Phi|\Phi'$ we assign
the functor $\underline\Gamma_{Z\setminus Z'}$, and to
any pair $(\underline Z_1,\underline Z_2)$ of elements
of $\Phi|\Phi'$ with $\underline Z_1\leq\underline Z_2$,
we attach the composition
\set\begin{equation}\label{eq_compose-gammas}
\underline\Gamma_{Z_1\setminus Z'_1}
\xrightarrow{\ \gamma_{Z_1\setminus Z'_1,Z_2\setminus Z'_1}\ }
\underline\Gamma_{Z_2\setminus Z'_1}
\xrightarrow{\ \gamma_{Z_2\setminus Z'_2,Z_2\setminus Z'_1}\ }
\underline\Gamma_{Z_2\setminus Z'_2}
\end{equation}
(notation of \eqref{subsec_variant-i-j}). We then define
the functor
$$
\underline\Gamma_{\Phi|\Phi'}:=\colim_{(Z,Z')\in\Phi|\Phi'}
\underline\Gamma_{Z\setminus Z'}:\cA\Mod\to\cA\Mod
$$
where the transition maps in the colimit are the
maps \eqref{eq_compose-gammas}, whence a derived
functor
$$
R\underline\Gamma_{\Phi|\Phi'}:\sD^+(\cA\Mod)\to\sD^+(\cA\Mod)
$$
and we have a natural isomorphism of functors on $\cA\Mod$ :
$$
R^i\underline\Gamma_{\Phi|\Phi'}\isom
\colim_{(Z,Z')\in\Phi|\Phi'}R^i\underline\Gamma_{Z\setminus Z'}
\qquad
\text{for every $i\in\N$}.
$$

\sset\subsubsection{}
Keep the notation of \eqref{subsec_family-of-supps}.
As a special case, we may take $\Phi':=\{\emptyset\}$,
and we let
$$
\underline\Gamma_\Phi:=\underline\Gamma_{\Phi|\{\emptyset\}}:
\cA\Mod\to\cA\Mod.
$$
If $\Phi$ and $\Phi'$ are arbitrary family of
supports, we have maps of partially ordered sets
$$
\Phi\xleftarrow{}\Phi|\Phi'\to\Phi'
\quad :\quad
Z\mapsfrom(Z,Z')\mapsto Z'.
$$
Suppose moreover that $\Phi'\subset\Phi$; then
both of these maps are cofinal, by lemma
\ref{lem_filtered-final}(i), and with lemma
\ref{subsec_formalities} and proposition
\ref{prop_derive-undergamma}(ii), we deduce a
short exact sequence of complexes
\set\begin{equation}\label{eq_use-this}
0\to\underline\Gamma_{\Phi'}\cI\to\underline\Gamma_\Phi\cI
\to\underline\Gamma_{\Phi|\Phi'}\cI\to 0
\end{equation}
for every flabby $\cA$-module $\cI$, whence a
distinguished triangle of functors on $\sD^+(\cA\Mod)$ :
$$
R\underline\Gamma_{\Phi'}\to R\underline\Gamma_\Phi
\to R\underline\Gamma_{\Phi|\Phi'}\to
R\underline\Gamma_{\Phi'}[1].
$$

\sset\subsubsection{}\label{subsec_supp-and-cod}
Keep the notation of \eqref{subsec_family-of-supps}.
Suppose now that $X$ is sober and locally noetherian, and
moreover that we are given a {\em weak codimension function}
on $X$, {\em i.e.} a map
$$
c:X\to\Z
$$
such that $c(x)>c(y)$ whenever $x\neq y$ and $x$ is
a specialization of $y$ in $X$. To the weak codimension
function $c$ we attach a system $(\Phi^p~|~p\in\N)$
of families of supports, by the rule :
$$
\Phi^p:=\{Z\subset X~|~
\text{$Z$ is closed in $X$, and $c(z)\geq p$ for every $z\in Z$}\}
\qquad
\text{for every $p\in\N$}.
$$
For every $x\in X$, let also $j_x:X(x)\to X$ and
$i_x:\{x\}\to X$ be the inclusion maps (notation
of definition \ref{def_special}(iii)). Moreover,
let us set $\Phi^p(x):=\{Z\cap X(x)~|~Z\in\Phi^p\}$,
for every $p\in\Z$.

\begin{lemma}\label{lem_local-situation}
In the situation of \eqref{subsec_supp-and-cod},
the following holds for every $x\in X$ :
\begin{enumerate}
\item
The topological space $X(x)$ is spectral, noetherian
and finite dimensional, and the function $c$ is bounded
on $X(x)$.
\item
There exists a natural isomorphism of functors on
$\sD^+(j^*_x\cA\Mod)$ :
$$
j^*_xR\underline\Gamma_{\Phi^p}\isom R\underline\Gamma_{\Phi^p(x)}
\qquad
\text{for every $p\in\Z$}.
$$
\item
$(\Phi^p(x)~|~p\in\Z)$ is the system of families of
supports on $X(x)$ attached to the weak codimension
function $c_{|X(x)}$.
\end{enumerate}
\end{lemma}
\begin{proof}(i): Since the topology of $X(x)$ is
induced from that of $X$, it is clear that $X(x)$
is noetherian, and we have already seen that it
is spectral (remark \ref{rem_specialize}(i)).
Especially, $X(x)$ admits finitely many maximal
points (remark \ref{rem_specialize}(iii)); if
$\eta_1,\dots,\eta_k$ are these maximal points,
set $C:=\max(c(\eta_1),\dots,c(\eta_k))$.
Then it is clear that the restriction to $X(x)$
of the function $c$ takes its values in the
range $C,C+1,\dots,c(x)$, and it follows easily
that the length of any descending sequence of
points in $(X(x),\leq)$ is bounded by $c(x)-C$.

(ii) follows easily from proposition
\ref{prop_derive-undergamma}(v.c), since every closed
subset of a locally noetherian topological space is
constructible.

(iii): If $Z$ is any closed subset of $X(x)$, denote
by $\bar Z$ the topological closure of $Z$ in $X$;
since $X(x)$ is a pro-constructible subset of $X$,
the same holds for $Z=X(x)\cap\bar Z$. Then $\bar Z$
is the set of specializations in $X$ of the points
of $Z$ (proposition \ref{prop_closed-under-spec}(i)).
Especially, if $c(z)\geq p$ for every $z\in Z$, then
$c(z)\geq p$ for every $z\in\bar Z$, and then clearly
$Z\in\Phi^p(x)$. The assertion follows immediately.
\end{proof}

\begin{proposition}\label{prop_supp-and-cod}
In the situation of \eqref{subsec_supp-and-cod},
we have a natural isomorphism of functors
$$
R^i\underline\Gamma_{\Phi^p|\Phi^{p+1}}\isom
\bigoplus_{c(x)=p}j_{x*}\circ R^i\underline\Gamma_{\{x\}}\circ j^*_x
\qquad
\text{for every $p\in\Z$ and $i\in\N$}.
$$
\end{proposition}
\begin{proof} Since the assertion is local on $X$,
we may assume that $X$ is noetherian. We consider
first the case where $i=0$, and we fix $p\in\Z$.
Let $Z\in\Phi^p$ be any element; recall that $Z$
has a finite number of irreducible components
$Z_1,\dots,Z_k$ (proposition \ref{prop_noetherian}).
Let $\eta_i$ be the generic point of $Z_i$, so that
$c(\eta_i)\geq p$ for every $i=1,\dots,k$.
For every $i,j=1,\dots,k$ with $i\neq j$,
set $Z_{ij}:=Z_i\cap Z_j$; then, for every such
$i,j$, every irreducible component of $T$ of
$Z_{ij}$ and every $t\in T$, we must have
$c(t)\geq c(\eta_T)>c(\eta_i)$, if $\eta_T$ denotes
the generic point of $T$. We conclude that
$Z_{ij}\in\Phi^{p+1}$ for every $i,j=1,\dots k$ with
$i\neq j$. Set also
$$
Z^{>p}:=\bigcup_{c(\eta_i)>p}Z_i\cup\bigcup_{i\neq j}Z_{ij}
\qquad\text{and}\qquad
\Psi^{p+1}:=\{(Z,Z')\in\Phi^p|\Phi^{p+1}~|~Z^{>p}\subset Z'\subset Z\}.
$$
Since the subset $\Psi^p$ is clearly cofinal in
$\Phi^p|\Phi^{p+1}$, we deduce a natural isomorphism
$$
\underline\Gamma_{\Phi^p|\Phi^{p+1}}\isom
\colim_{(Z,Z')\in\Psi^p}\underline\Gamma_{Z\setminus Z'}.
$$
However, say that $(Z,Z')\in\Psi^p$; we notice that
the closed subset $V:=Z\!\setminus\!Z'$ of $X\!\setminus\!Z'$
is the disjoint union of its irreducible components
$V_1,\dots,V_t$, hence the topology induced by $X$
on $V$ is the coproduct ({\em i.e.} disjoint union)
of the topologies induced by $X$ on $V_1,\dots,V_t$,
and there follows a natural isomorphism
\set\begin{equation}\label{eq_decompose-irred}
\underline\Gamma_V\isom
\bigoplus_{i=1}^t\underline\Gamma_{V_i}.
\end{equation}
Moreover, say that $\underline Z_1,\underline Z_2\in\Psi^p$
and $\underline Z_1\leq\underline Z_2$, and denote by
$V_1,\dots,V_s$ (resp. $W_1,\dots,W_t$) the irreducible
components of $V:=Z_1\!\setminus\!Z'_1$ (resp. of
$W:=Z_2\!\setminus\!Z'_2$).
A simple inspection of the definitions reveals that $s\leq t$
and -- after reordering of the $W_i$ -- the subset $W_i$ is
open in $V_i$ for every $i=1,\dots,s$. Moreover, under the
isomorphism \eqref{eq_decompose-irred}, the transition
map $\underline\Gamma_V\to\underline\Gamma_W$ is identified
with the composition
$$
\bigoplus_{i=1}^s\underline\Gamma_{V_i}
\xrightarrow{\ \bigoplus_{i=1}^s\gamma_{V_i,W_i}\ }
\bigoplus_{i=1}^s\underline\Gamma_{W_i}\to
\bigoplus_{i=1}^t\underline\Gamma_{W_i}
$$
(where the second map is the natural monomorphism).
Hence, let $\Delta^p$ be the set of all locally
closed subsets $V$ of $X$ fulfilling the following
condition. There exists $x\in X$ with $c(x)=p$, and
$V$ is a non-empty open subset of the topological
closure $\bar{\{x\}}$ of $\{x\}$ in $X$. We endow
$\Delta^p$ with a partial ordering, by declaring
that, for any $V,V'\in\Delta^p$ we have $V\leq V'$
if $V'\subset V$ (notice that $\Delta^p$ is usually
not a filtered partially ordered set); from the
foregoing, we get a natural isomorphism
$$
\underline\Gamma_{\Phi|\Phi}\isom
\colim_{V\in\Delta^p}\underline\Gamma_V
$$
where, for every $V,V'\in\Delta^p$ such that $V\leq V'$,
the transition map
$\underline\Gamma_V\to\underline\Gamma_{V'}$ is again
given by the transformation $\gamma_{V,V'}$ corresponding
to the open immersion $V'\to V$.

Next, for every $x\in X$ such that $c(x)=p$, set
$\Delta_x:=\{V\in\Delta^p~|~x\in V\}$, and notice
that the partially ordered set $\Delta^p$ is the
coproduct ({\em i.e.} the disjoint union) of its
partially ordered subsets $\Delta_x$. We deduce
a natural isomorphism
$$
\underline\Gamma_{\Phi|\Phi}\isom\bigoplus_{c(x)=p}
\colim_{V\in\Delta_x}\underline\Gamma_V.
$$
Thus, we come down to checking the following :

\begin{claim} For every $p\in\Z$ and every $x\in X$
such that $c(x)=p$, there exists a natural isomorphism
of functors
$$
\colim_{V\in\Delta^p_x}\underline\Gamma_V\isom
j_{x*}\circ\underline\Gamma_{\{x\}}\circ j^*_x.
$$
\end{claim}
\begin{pfclaim} Let $\Sigma_x$ be the set of all
quasi-compact open neighborhoods of $x$ in $X$, and
endow $\Sigma_x$ with the partial ordering such that,
for every $U,U'\in\Sigma_x$ we have $U\leq U'$ if
$U'\subset U$; then $\Sigma_x$ is filtered, and there
follows a morphism of partially ordered sets
$$
\Sigma_x\to\Delta_x
\qquad
U\mapsto U\cap\bar{\{x\}}
$$
and using lemma \ref{lem_filtered-final}(i), it is
easily seen that this map is cofinal. Now, for every
$U\in\Sigma_x$, let $j_U:U\to X$ be the open immersion;
we have a natural identification
$$
\underline\Gamma_{U\cap\bar{\{x\}}}\isom
j_{U*}\circ j_U^*\circ\underline\Gamma_{\bar{\{x\}}}
$$
whence a natural isomorphism
$$
\colim_{V\in\Delta_x}\underline\Gamma_V=
\colim_{U\in\Sigma_x}
j_{U*}\circ j_U^*\circ\underline\Gamma_{\bar{\{x\}}}.
$$
However, $X(x)$ is the limit of the filtered system
of spectral spaces $\Sigma_x$, so the claim follows
easily from proposition \ref{prop_dir-im-and-colim}(ii)
(details left to the reader).
\end{pfclaim}

Next, let $i\in\N$ be arbitrary, let $\cF$ be any
$\cA$-module, and pick any resolution $\cF\to\cI^\bullet$
consisting of injective $\cA$-modules; by the foregoing
case, we get a natural isomorphism of $\cA$-modules
$$
R^i\underline\Gamma_{\Phi^p|\Phi^{p+1}}\cF\isom
\bigoplus_{c(x)=p}
H^i(j_{x*}\circ\underline\Gamma_{\{x\}}\circ j^*_x\cI^\bullet).
$$
However, on the one hand, $j^*\cI^\bullet$ is a
complex of qc-flabby $j^*\cA$-modules (claim
\ref{cl_old-flame}); on the other hand, qc-flabby
abelian sheaves are acyclic for the functor
$\underline\Gamma_{\{x\}}$ (lemmata
\ref{lem_flabby-Gamma_Z}(iii.b) and
\ref{lem_charact-coh-sob-sp}(vi)), so we get a
natural isomorphism
$$
H^i(j_{x*}\circ\underline\Gamma_{\{x\}}\circ j^*_x\cI^\bullet)
\isom
H^i(j_{x*}\circ R\underline\Gamma_{\{x\}}\circ j^*_x\cF)
\qquad
\text{for every $i\in\N$ and every $x\in X$}.
$$
Lastly, notice that $R\underline\Gamma_{\{x\}}\circ j^*_x\cF$
is a complex of $j^*\cA$-modules supported on $\{x\}$,
and the functor $j_{x*}$ is exact on the full subcategory
of $j^*\cA\Mod$ consisting of such sheaves; thus, we
have a natural isomorphism
$$
H^i(j_{x*}\circ R\underline\Gamma_{\{x\}}\circ j^*_x\cF)
\isom
j_{x*}\circ R^i\underline\Gamma_{\{x\}}\circ j^*_x\cF
\qquad
\text{for every $i\in\N$ and every $x\in X$}.
$$
Summing up, we obtain the proposition.
\end{proof}

\begin{definition}\label{def_Cousin}
Let $(X,\cA)$ be a ringed topological space, with $X$
sober and locally noetherian, and $c$ a weak codimension
function on $X$.

(i)\ \
We say that an object $K^\bullet$ of $\sD^+(\cA\Mod)$
is a {\em $c$-Cohen-Macaulay complex}, if we have
$$
R^i\underline\Gamma_{\{x\}}K^\bullet_{|X(x)}=0
\qquad
\text{for every $x\in X$ and every integer $i\neq c(x)$}.
$$
We denote by
$$
\cA\CM^c
$$
the full subcategory of $\sD^+(\cA\Mod)$ whose objects
are the $c$-Cohen-Macaulay complexes.

(ii)\ \
We say that an object $C^\bullet$ of $\sC^+(\cA\Mod)$
is a {\em Cousin complex}, if there exists a family
$(M_x~|~x\in X)$ where $M_x$ is an $i_x^*\cA$-module
for every $x\in X$ (notation of \eqref{subsec_supp-and-cod}),
such that
$$
C^p=\bigoplus_{c(x)=p}i_{x*}M_x
\qquad
\text{for every $p\in\Z$}.
$$
The full subcategory of $\sC^+(\cA\Mod)$ whose objects
are the Cousin complexes is denoted
$$
\cA\text{-}\mathbf{Cousin}.
$$
\end{definition}

\begin{lemma}\label{lem_forget-cousin}
In the situation of definition {\em\ref{def_Cousin}},
we have :
\begin{enumerate}
\item
$R\underline\Gamma_{\{x\}}\circ j^*_x\circ i_{y*}=0$
for every $x,y\in X$ such that $x\neq y$.
\item
If $C^\bullet$ is a Cousin complex of $\cA$-modules,
then $C^p$ is qc-flabby, for every $p\in\Z$.
\end{enumerate}
\end{lemma}
\begin{proof}(i): Notice that any $i^*_y\cA$-module $M$
is trivially qc-flabby, so the same holds for $i_{y*}M$
(lemma \ref{lem_franziska}(iv)), and therefore also for
$j^*_xj_{y*}M$ (claim \ref{cl_old-flame}). Since qc-flabby
$\cA$-modules are acyclic for $\underline\Gamma_{\{x\}}$
(lemma \ref{lem_flabby-Gamma_Z}(iii.a)), we conclude
that $R\underline\Gamma_{\{x\}}j^*_xi_{y*}$ is the derived
functor of $\underline\Gamma_{\{x\}}j^*_xi_{y*}$; but
it is easily seen that the latter vanishes, whence the
assertion.

(ii): We have just seen that $i_*M$ is qc-flabby, for
every $i^*_x\cA$-module $M$; on the other hand, all
direct sums of qc-flabby $\cA$-modules are qc-flabby
(claim \ref{cl_old-flame}), whence the contention.
\end{proof}

\sset\subsubsection{}\label{subsec_forget-cousin}
As $j^*_x$ and $R^i\underline\Gamma_{\{x\}}$ commute with
all direct sums (proposition
\ref{prop_derive-undergamma}(v.b)), lemma
\ref{lem_forget-cousin}(i) implies that the localization
$\sC^+(\cA\Mod)\to\sD^+(\cA\Mod)$ restricts to a functor
$$
F:\cA\text{-}\mathbf{Cousin}\to\cA\CM^c.
$$

\begin{theorem}\label{th_Cousin}
The functor $F$ of \eqref{subsec_forget-cousin}
is an equivalence of categories.
\end{theorem}
\begin{proof} We construct a quasi-inverse to $F$, as
follows. First, we endow every object $C^\bullet$ of
$\sC(\cA\Mod)$, with a descending filtration
$\Fil^\bullet_\Phi C^\bullet$, by the rule
$$
\Fil_\Phi^pC^\bullet:=\underline\Gamma_{\Phi^p}C^\bullet
$$
where $(\Phi^p~|~p\in\Z)$ is the system of families of
supports associated with the weak codimension function $c$,
as in \eqref{subsec_supp-and-cod}. We also denote by
$\gr^\bullet_\Phi C^\bullet$ the resulting associated
graded object of $\sC(\cA\Mod)$. 
Next, for every object $K^\bullet$ of $\sD^+(\cA\Mod)$
we choose a quasi-isomorphism
$\beta^\bullet_K:K^\bullet\isom\cI^\bullet_K$, with
$\cI^\bullet_K$ a bounded below complex of injective
$\cA$-modules. Following \eqref{subsec_spec-seq-fil-cpx},
we denote
$E(\cI_K)_\bullet^{\bullet\bullet}$ the spectral sequence
arising from the filtered complex
$(\cI^\bullet_K,\Fil^\bullet_\Phi\cI_K^\bullet)$. In light of
\eqref{eq_use-this} we have a natural identification
$$
\gr^p_\Phi\cI^\bullet_K\isom
\underline\Gamma_{\Phi^p|\Phi^{p+1}}\cI^\bullet_K
$$
and combining with remark \ref{rem_low-terms}(ii),
we get a natural isomorphism
$$
E(\cI_K)_1^{pq}\isom
R^{p+q}\underline\Gamma_{\Phi^p|\Phi^{p+1}}K^\bullet.
$$
Taking into account proposition \ref{prop_supp-and-cod}
we deduce that
$$
G(K)^\bullet:=(E(\cI_K)_1^{\bullet 0},d(\cI_K)_1^{\bullet 0})
$$
is a Cousin complex, for every $q\in\Z$. Moreover,
every morphism $\phi^\bullet:K^\bullet\to L^\bullet$ in
the category
$\sD^+(\cA\Mod)$ can be represented by a unique morphism
$\cI^\bullet_\phi:\cI_K^\bullet\to\cI_L^\bullet$ in
$\Hot(\cA\Mod)$
(notation of definition \ref{def_homotopies}(ii)), and
notice that $\cI^\bullet_\phi$ is a morphism of filtered
complexes; in light of \eqref{subsec_hot-and-SpSeq}
we deduce a well defined morphism
$G(\phi)^\bullet:G(K)^\bullet\to G(L)^\bullet$. In this way,
we obtain a functor
$$
G:\sD^+(\cA\Mod)\to\cA\text{-}\mathbf{Cousin}
$$
and we claim that the restriction of $G$ to the subcategory
$\cA\CM^c$ is the sought quasi-inverse.

We check first that $G\circ F$ is naturally isomorphic
to the identity endofunctor of $\cA\text{-}\mathbf{Cousin}$.
To this aim, let us remark :

\begin{claim}\label{cl_self-aggrandize}
Let $K^\bullet$ be any $c$-Cohen-Macaulay complex. We have :
\begin{enumerate}
\item
$E(\cI_K)_1^{pq}=0$ for every $p,q\in\Z$ with $q\neq 0$.
\item
If $K^\bullet=FC^\bullet$ for a Cousin complex $C^\bullet$,
the map $\beta^\bullet_K$ induces a natural isomorphism
$$
\gr_\Phi^pC^\bullet\isom
E(\cI_K)_1^{p0}[-p]\isom C^p[-p]
\qquad
\text{for every $p\in\Z$}.
$$
\end{enumerate}
\end{claim}
\begin{pfclaim}(i): Indeed, from proposition
\ref{prop_supp-and-cod} we get a natural isomorphism
\set\begin{equation}\label{eq_self-import}
H^{p+q}(\gr^p_\Phi\cI^\bullet_C)\isom\bigoplus_{c(x)=p}
j_{x*}\circ R^{p+q}\underline\Gamma_{\{x\}}K^\bullet_{|X(x)}
\end{equation}
which vanishes for $q\neq 0$, since $K^\bullet$ is
$c$-Cohen-Macaulay.

(ii): Let $C^\bullet$ and $(M_x~|~x\in X)$ be as in
definition \ref{def_Cousin}(ii); from \eqref{eq_use-this}
and lemma \ref{lem_forget-cousin}(ii), it follows that
the natural map
$$
\gr^\bullet_\Phi C^\bullet\to
\underline\Gamma_{\Phi^p|\Phi^{p+1}}C^\bullet
$$
is an isomorphism in $\sC^+(\cA\Mod)$.
Also, combining with claim \ref{cl_old-flame},
we see that $C^\bullet_{|X(x)}$ is a complex of qc-flabby
$j^*_x\cA$-modules. By lemma \ref{lem_flabby-Gamma_Z}(iii.a)
we deduce that the natural map
$\underline\Gamma_{\{x\}}C^\bullet_{|X(x)}\to
R\underline\Gamma_{\{x\}}C^\bullet_{|X(x)}$ is an
isomorphism in $\sD^+(\cA\Mod)$, for every $x\in X$.
Moreover, since the functor $\underline\Gamma_{\{x\}}$
commutes with arbitrary direct sums (lemma
\ref{lem_flabby-Gamma_Z}(iii.b)), lemma
\ref{lem_forget-cousin}(i) says that
$$
\underline\Gamma_{\{x\}}C^\bullet_{|X(x)}=j^*_xi_{x*}M_x[-c(x)]
\qquad
\text{for every $x\in X$}.
$$
Taking into account \eqref{eq_self-import} (and again,
proposition \ref{prop_supp-and-cod}), the assertion follows.
\end{pfclaim}

Now, let $C^\bullet$ be a cousin complex, and set
$D^\bullet:=GF(C)^\bullet$. From claim \ref{cl_self-aggrandize}
we obtain a natural isomorphism
$$
\gr^p_\Phi\cI^\bullet_C\isom C^p[-p]
\qquad
\text{for every $p\in\Z$}.
$$
We need to show that, under this identification, the
differential $d^p_D$ of $D^p:=H^p(\gr^p\cI^\bullet_C)$
corresponds to the differential $d^p_C$ of $C^\bullet$.
However, on the one hand we have a commutative ladder
$$
\xymatrix{
0 \ar[r] & \gr^{p+1}_\Phi C^\bullet \ar[r] \ar[d] &
\Fil^p_\Phi C^\bullet/\Fil^{p+2}_\Phi C^\bullet \ar[r] \ar[d] &
\gr^p_\Phi C^\bullet \ar[r] \ar[d] & 0 \\
0 \ar[r] & \gr^{p+1}_\Phi\cI_C^\bullet \ar[r] &
\Fil^p_\Phi\cI_C^\bullet/\Fil^{p+2}_\Phi\cI_C^\bullet \ar[r] &
\gr^p_\Phi\cI_C^\bullet \ar[r] & 0
}$$
and claim \ref{cl_self-aggrandize}(ii) says that both the
first and third vertical arrows are isomorphisms, so the
top row of is isomorphic to the bottom row. Moreover, we
see that the top row is the natural short exact sequence
$$
0\to C^{p+1}[-p-1]\to t^{\leq p+1}t^{\geq p}C^\bullet\to
C^p[-p]\to 0
$$
where $t^{\leq p+1}$ and $t^{\geq p}$ are the brutal
truncation functors of \eqref{subsec_brutal-truncate}.
Then it is clear that the boundary map
$H^p(C^p[-p])\to H^{p+1}(C^{p+1}[-p-1])$ in degree
$p$ associated with this short exact sequence, is none
else than the differential $d_C^p$. Now the contention
follows from the naturality of the boundary map, together
with remark \ref{rem_low-terms}(ii).

Next, let $K^\bullet$ be any $c$-Cohen-Macaulay complex,
and set
$$
(\cJ_K,\Fil^\bullet\cJ_K):=
\Dec(\cI_K,\Fil^\bullet_\Phi\cI_K^\bullet)
$$
(notation of \eqref{subsec_decalage}). By
proposition \ref{prop_decalage}(i,ii) and claim
\ref{cl_self-aggrandize}(i), we have a natural
epimorphism
$$
u^\bullet:E(\cJ_K)_0^\bullet\to E(\cI_K)_1^\bullet
\qquad
\text{in $\sC^+(\cA\Mod)$}
$$
which induces an isomorphism in $\sD^+(\cA\Mod)$.
Moreover, from claim \ref{cl_self-aggrandize}(i)
we see that the subcomplex $(E(\cJ_K)_0^{p,n-p}~|~n\in\Z)$
lies in $\Ker\,u^\bullet$, for every $p\neq 0$,
so $u^\bullet$ restricts to an isomorphism
\set\begin{equation}\label{eq_finish-cousin}
E(\cJ_K)_0^{0,\bullet}=\gr^0\cJ_K^\bullet\isom G(K)^\bullet
\qquad
\text{in $\sD^+(\cA\Mod)$}
\end{equation}
(see remark \ref{rem_low-terms}(i)). We point out :

\begin{claim}\label{cl_fil-acyclic}
$\Fil^1\cJ_K^\bullet$ and $\cJ_K^\bullet/\Fil^0\cJ_K^\bullet$
are acyclic.
\end{claim}
\begin{pfclaim} The assertion can be checked on the
stalks over the points of $X$; taking into account
lemma \ref{lem_local-situation}(ii,iii), we may then
replace $X$ by $X(x)$, for any point $x\in X$, and
assume from start that $c$ is bounded on $X$
(lemma \ref{lem_local-situation}(i)). In this case,
it is clear that the filtration
$\Fil^\bullet_\Phi\cI^\bullet_K$ is finite, for every
complex $K^\bullet$ of $\cA$-modules, and then the
same holds for the filtration $\Fil^\bullet\cJ^\bullet_K$
(see \eqref{subsec_decalage}). In this situation,
a simple induction argument reduces to checking that
the associated subquotients $\gr^p\cJ_K^\bullet$ are
acyclic, for every $p\neq 0$. However, combining
remarks \ref{rem_vanish-reproduces}(i) and
\ref{rem_low-terms}(ii), proposition \ref{prop_decalage}(ii)
and claim \ref{cl_self-aggrandize}(i) we may compute
$$
H^q\gr^p\cJ_K^\bullet\isom E(\cJ_K)_1^{p,q-p}\isom
E(\cI_K)_2^{p+q,-p}=0
\qquad
\text{for every $p,q\in\Z$ with $p\neq 0$}
$$
whence the claim.
\end{pfclaim}

Claim \ref{cl_fil-acyclic} yields a natural isomorphism
$$
\gr^0\cJ^\bullet_K\isom\cJ^\bullet_K=\cI^\bullet_K
\qquad
\text{in $\sD^+(\cA\Mod)$}
$$
and combining with \eqref{eq_finish-cousin}, we
finally get an isomorphism
$$
K^\bullet\isom FG(K)^\bullet
\qquad
\text{in $\sD^+(\cA\Mod)$}.
$$
Now the theorem follows from proposition
\ref{prop_fullfaith-adjts}(i).
\end{proof}

\subsection{Duality over coherent schemes}\label{subsec_duals}
If $X$ is a locally coherent scheme (definition
\ref{def_coh-schemes}(i)) and $\cF$ is an $\cO_{\!X}$-module,
we define the {\em dual $\cO_{\!X}$-module}
$$
\cF^\vee:=\cHom_{\cO_{\!X}}(\cF,\cO_{\!X}).
$$
Moreover, for every morphism of $\cO_{\!X}$-modules $\phi:\cF\to\cG$,
we denote by $\phi^\vee:\cG^\vee\to\cF^\vee$ the induced (transpose)
morphism. As usual, there is a natural morphism of $\cO_{\!X}$-modules:
$$
\beta_\cF:\cF\to\cF^{\vee\vee}.
$$

\begin{definition}\label{def_reflex}
Let $X$ be a locally coherent scheme, and $\cF$ an
$\cO_{\!X}$-module.
\begin{enumerate}
\item
We say that $\cF$ is {\em reflexive at a point\/} $x\in X$ if
there exists an open neighborhood $U\subset X$ of $x$ such
that $\cF_{|U}$ is a coherent $\cO_{\!U}$-module, and
$\beta_{\cF|U}$ is an isomorphism.
\item
We say that $\cF$ is {\em reflexive\/} if it is reflexive
at all points of $X$.
\item
We denote by $\cO_{\!X}\bRflx$ the full subcategory of the category
$\cO_{\!X}\Mod$, consisting of all the reflexive $\cO_{\!X}$-modules.
It contains $\cO_{\!X}\Mod_\mathrm{lfft}$ as a full subcategory
(see \eqref{sec_various-O-mod}).
\item
Suppose that $X=\Spec\,R$ for a coherent ring $R$, and
$\cF=M^\sim$ for an $R$-module $M$. Then we say that $M$
is a {\em reflexive} $R$-module if $\cF$ is a reflexive
$\cO_{\!X}$-module.
\end{enumerate}
\end{definition}

\begin{lemma}\label{lem_refl-is-local}
Let $X$ be a locally coherent scheme, $x\in X$ any point,
and $\cF$ a coherent $\cO_{\!X}$-module. We have :
\begin{enumerate}
\item
$\cF^\vee$ is a coherent $\cO_{\!X}$-module.
\item
The following conditions are equivalent:
\begin{enumerate}
\item
$\cF$ is reflexive at the point $x$.
\item
$\cF_{\!x}$ is a reflexive $\cO_{\!X,x}$-module.
\item
The map $\beta_{\cF,x}:\cF_{\!x}\to(\cF^{\vee\vee})_x$
is an isomorphism.
\romanenu
\end{enumerate}
\end{enumerate}
\end{lemma}
\begin{proof}(i): The assertion is local on $X$,
hence we may assume that there exists a presentation
$\cO_{\!X}^{\oplus n}\to\cO_{\!X}^{\oplus m}\to\cF\to 0$,
whence -- after taking duals -- a left exact sequence
$0\to\cF^\vee\to\cO_{\!X}^{\oplus m}\to\cO_{\!X}^{\oplus n}$.
Since $\cO_{\!X}$ is coherent, the assertion follows.

(ii) shall be left to the reader.
\end{proof}

\begin{lemma}\label{lem_dual-is-rflx}
Let $X$ be a reduced locally coherent scheme, and
$\cF$ a quasi-coherent $\cO_{\!X}$-module of finite
type. Endow the set $\Max\,X$ of maximal points of $X$
with the topology induced by $X$, and consider the
following conditions :
\begin{enumerate}
\alphaenu
\item
There exists a quasi-compact open immersion
$j:U\to X$ with dense image, such that $j^*\cF$ is a
locally free $\cO_{\!U}$-module.
\item
There exists a quasi-compact open immersion
$j:U\to X$ with dense image, such that $j^*\cF$ is a
coherent $\cO_{\!U}$-module.
\item
The rank function of $\cF$ :
$$
X\to\N
\qquad
x\mapsto\dim_{\kappa(x)}\kappa(x)\otimes_{\cO_{\!X,x}}\cF_x
$$
restricts to a locally constant function on $\Max\,X$.
\item
$\cF^\vee$ is a reflexive $\cO_{\!X}$-module and
$\beta^\vee_{\cF}\circ\beta_{\cF^\vee}=\one_{\cF^\vee}$.
\end{enumerate}
Then {\em
(a)$\Rightarrow$(b)$\Rightarrow$(c)$\Rightarrow$(d)} and
if $X$ is coherent we have {\em (c)$\Rightarrow$(a)} as well.
\end{lemma}
\begin{proof} Let us check first that (c)$\Rightarrow$(a)
in case $X$ is coherent. Indeed, notice that $j$ is
quasi-compact if and only if $U$ is retro-compact in
$X$, and the latter holds if and only if $U$ is open
and constructible in $X$ (lemma
\ref{lem_sorite-construct}(iv.b,v.c)). It follows that
the existence of $j$ with the sought properties can
be checked locally on $X$, hence we may assume that
$X$ is affine. By proposition
\ref{prop_coher-then-noether}(i) we may find a coherent
$\cO_{\!X}$-module $\cG$ and an $\cO_{\!X}$-linear
epimorphism $\phi:\cG\to\cF$. Now, fix any $\eta\in\Max\,X$
and write $\cK:=\Ker\,\phi$ as the filtered union of the
system $(\cK_i~|~i\in I)$ of its quasi-coherent
$\cO_{\!X}$-submodules of finite type; since
$\cO_{\!X,\eta}$ is a field, there exists $i\in I$
such that $\cK_{i,\eta}=\cK_\eta$. Thus, after
replacing $\cG$ by $\cG/\cK_i$, we may assume that
$\phi_\eta$ is an isomorphism. Moreover, by
\cite[Th.4.10]{Mat} we may find an affine open
neighborhood $U_\eta$ of $\eta$ in $X$, an integer
$r\in\N$ and an isomorphism
$\cG_{|U_\eta}\isom\cO^{\oplus r}_{\!U_\eta}$, and by assumption
(c) we may also assume that the rank function of $\cF$
is constant on $\Max\,U_\eta$. In this case, $\phi_\tau$
is an isomorphism for every $\tau\in\Max\,U_\eta$, so
$\cK_\tau=0$ for every such point $\tau$. Since $X$ is
reduced and $\cG_{|U_\eta}$ is a locally free
$\cO_{\!U_\eta}$-module, it follows that $\cK_{|U_\eta}=0$.
By proposition \ref{prop_max-is-qc}, we may then find a
finite subset $\Sigma\subset\Max\,X$ such that
$\bigcup_{\eta\in\Sigma}\Max\,U_\eta=\Max\,X$, and it
is easily seen that the open subset
$U:=\bigcup_{\eta\in\Sigma}U_\eta$ fulfills condition (a).

(a)$\Rightarrow$(b) is obvious.

(b)$\Rightarrow$(c): Notice first that any $U$ as in
(b) contains $\Max\,U$, by proposition
\ref{prop_closed-under-spec}(ii). Then, since $X$ is reduced,
\cite[Th.4.10]{Mat} implies that, for every $\eta\in\Max\,X$,
we may find an affine open neighborhood $U_\eta$ of $\eta$
in $X$ such that $\cF_{|U_\eta}$ is a locally free
$\cO_{\!U_\eta}$-module, whence (c).

(c)$\Rightarrow$(d): Let us show first that $\cF^\vee$
is coherent, under condition (c). The assertion is
local on $X$, hence we may assume that $X$ is affine,
in which case, as we have already seen, (c) implies
that there exists a quasi-compact open immersion
$j:U\to X$ fulfilling condition (a), and we also
know that any such $U$ contains $\Max\,X$. 
By proposition \ref{prop_coher-then-noether}(i),
we may then find a coherent $\cO_{\!X}$-module $\cG$
and an $\cO_{\!X}$-linear epimorphism $\phi:\cG\to\cF$
such that $j^*\phi$ is an isomorphism. There follows
a left exact sequence
$0\to\cF^\vee\to\cG^\vee\to(\Ker\,\phi)^\vee$, and
$\cG^\vee$ is coherent, by lemma \ref{lem_refl-is-local}(i).
Since $\Max\,X\subset U$, we have $(\Ker\,\phi)_\eta=0$
for every $\eta\in\Max\,X$, and since $\cO_{\!X}$ is reduced,
it follows easily that $(\Ker\,\phi)^\vee=0$, and the claim
follows. Now, since $X$ is reduced, the only associated points
of $\cO_{\!X}$ are the maximal points of $X$ ({\em i.e.}
$\cO_{\!X}$ has no imbedded points). It follows easily that,
for every quasi-coherent $\cO_{\!X}$-module $\cG$, the dual
$\cG^\vee$ satisfies condition $S_1$ (see definition
\ref{def_Ass}(ii)). Next, rather generally, let $\cM$ be any
$\cO_{\!X}$-module; directly from the definitions one derives
the identity:
\set\begin{equation}\label{eq_left-inverse}
\beta^\vee_{\!\cM}\circ\beta_{\!\cM^\vee}=\one_{\cM^\vee}.
\end{equation}
It remains therefore only to show that $\beta^\vee_\cF$ is
a right inverse for $\beta_{\cF^\vee}$ when $\cF^\vee$ is
coherent. Since $\cF^{\vee\vee\vee}$ satisfies condition
$S_1$, it suffices to check that, for every maximal point
$\xi$, the induced map on stalks
$$
\beta^\vee_{\cF,\xi}:\cF^{\vee\vee\vee}_\xi\to\cF^\vee_\xi
$$
is a right inverse for $\beta_{\cF^\vee,\xi}$. However, since
$\cO_{\!X,\xi}$ is a field, $\beta^\vee_{\cF,\xi}$ is a linear
map of $\cO_{\!X,\xi}$-vector spaces of the same finite
dimension, hence it is an isomorphism, in view of
\eqref{eq_left-inverse}.
\end{proof}

\begin{remark}\label{rem_co-repres}
The following observation is often useful. Suppose that $X$ is
a locally coherent scheme, $\cF$ an $\cO_{\!X}$-module, reflexive
at a given point $x\in X$. We can then choose a presentation
$\cO^{\oplus n}_{\!X,x}\to\cO^{\oplus m}_{\!X,x}\to\cF_x^\vee\to 0$,
and after dualizing we deduce a left exact sequence
\set\begin{equation}\label{eq_co-repres}
0\to\cF_x\to\cO^{\oplus m}_{\!X,x}\xrightarrow{u}\cO^{\oplus n}_{\!X,x}.
\end{equation}
Especially, if $\cO_{\!X,x}$ satisfies condition $S_1$ (in the sense
of definition \ref{def_Ass}(iii)), then the same holds for $\cF_x$.
For the converse, suppose additionally that $X$ is reduced, and let
$x\in X$ be a point for which there exists a left exact sequence
such as \eqref{eq_co-repres}; then lemma \ref{lem_dual-is-rflx}
says that $\cF$ is reflexive at $x$ : indeed,
$\cF_x\simeq(\Coker\,u^\vee)^\vee$.
\end{remark}

\begin{lemma}\label{lem_rflx-on-lim}
{\em (i)}\ \  
Let $f:X\to Y$ be a flat morphism of locally coherent schemes.
The induced functor $\cO_Y\Mod\to\cO_{\!X}\Mod$ restricts to a
functor
$$
f^*:\cO_Y\bRflx\to\cO_{\!X}\bRflx.
$$
\begin{enumerate}
\addenu
\item
Let $X_0$ be a  quasi-compact and quasi-separated scheme,
$(X_\lambda~|~\lambda\in\Lambda)$ a cofiltered family of
coherent $X_0$-schemes with flat transition morphisms
$\psi_{\lambda\mu}:X_\lambda\to X_\mu$ such that the
structure morphisms $X_\lambda\to X_0$ are affine, and set
$X:=\lim_{\lambda\in\Lambda}\,X_\lambda$. Then:
\begin{enumerate}
\item
$X$ is coherent.
\item
the natural functor:
$\Pscolim{\lambda\in\Lambda^o}\cO_{\!X_\lambda}\bRflx\to\cO_{\!X}\bRflx$
is an equivalence.
\end{enumerate}
\item
Suppose that $f$ is surjective, and let $\cF$ be any coherent
$\cO_Y$-module. Then $\cF$ is reflexive if and only if $f^*\cF$
is a reflexive $\cO_{\!X}$-module.
\end{enumerate}
\end{lemma}
\begin{proof} (i) follows easily from \cite[Lemma 2.4.29(i.a)]{Ga-Ra}.

(ii): For every $\lambda\in\Lambda$, denote by
$\psi_\lambda:X\to X_\lambda$ the natural morphism. Let $U\subset X$
be a quasi-compact open subset, $u:\cO_{\!U}^{\oplus n}\to\cO_{\!U}$
any morphism of $\cO_{\!U}$-modules; by corollary
\ref{cor_main-spectral}(ii.a) there exist
$\lambda\in\Lambda$ and a quasi-compact open subset
$U_\lambda\subset X_\lambda$ such that $U=\psi^{-1}(U_\lambda)$.
By \cite[Ch.IV, Th.8.5.2(i)]{EGAIV-3} we may then suppose that
$u$ descends to a homomorphism
$u_\lambda:\cO_{\!U_\lambda}^{\oplus n}\to\cO_{\!U_\lambda}$,
whose kernel is of finite type, since $X_\lambda$ is coherent.
Since the transition morphisms are flat, we have
$\Ker\,u=\psi^*_\lambda(\Ker\,u_\lambda)$, whence (ii.a).
Next, using \cite[Ch.IV, Th.8.5.2]{EGAIV-3} one sees easily
that the functor of (ii.b) is fully faithful and moreover,
every reflexive $\cO_{\!X}$-module $\cF$ descends to a coherent
$\cO_{\!X_\lambda}$-module $\cF_\lambda$ for some $\lambda\in\Lambda$.
For every $\mu\geq\lambda$ let
$\cF_{\!\mu}:=\psi_{\mu\lambda}^*\cF_\lambda$; since $\beta_{\cF}$
is an isomorphism, {\em loc.cit.} shows that $\beta_{\cF_\mu}$
is already an isomorphism for some $\mu\geq\lambda$, whence (ii.b).

(iii): By virtue of (i), we may assume that $f^*\cF$ is reflexive,
and we need to show that the same holds for $\cF$. However, the
natural map $f^*(\cF^{\vee\vee})\to(f^*\cF)^{\vee\vee}$ is an
isomorphism (proposition \ref{prop_replace-prop.12.3.5}(ii)),
hence $f^*\beta_\cF=\beta_{f^*\cF}$ is an isomorphism. Since $f$
is faithfully flat, we deduce that $\beta_\cF$ is an isomorphism,
as stated.
\end{proof}

\begin{proposition}\label{prop_extend-rflx}
Let $X$ be a coherent scheme, $U\subset X$ a quasi-compact
open subset. Then:
\begin{enumerate}
\item
If $X$ is reduced, the restriction functor
$$
\cO_{\!X}\bRflx\to\cO_{\!U}\bRflx
$$
is essentially surjective.
\item
Let 
$\cF$ be a reflexive $\cO_{\!X}$-module. If
$\delta'(x,\cO_{\!X})\geq 2$ for all $x\in X\setminus U$,
the natural map
$$
\cF\to j_*j^*\cF
$$
is an isomorphism. Especially, the restriction functor of\/
{\em (i)} is an equivalence.
\end{enumerate}
\end{proposition}
\begin{proof}(i): Given a reflexive $\cO_{\!U}$-module $\cF$,
lemma \ref{lem_extend-cohs}(ii) says that we can find a finitely
presented quasi-coherent $\cO_{\!X}$-module $\cG$ extending $\cF^\vee$;
since $X$ is coherent, $\cG$ is a coherent $\cO_{\!X}$-module, hence
the same holds for $\cG^\vee$, which extends $\cF$ and is reflexive
in light of lemma \ref{lem_dual-is-rflx}.

(ii): Notice that the first assertion easily implies
that the restriction functor of (i) is fully faithful,
so the second assertion follows from the first together
with (i).

Next, since the first assertion is local on $X$, we can
suppose that there exists a left exact sequence
$0\to\cF\to\cO^{\oplus m}_{\!X}\to\cO^{\oplus n}_{\!X}$
(see remark \ref{rem_co-repres}).
Since the functor $j_*$ is left exact, it then suffices
to prove the contention for the sheaves $\cO^{\oplus m}_{\!X}$
and $\cO^{\oplus n}_{\!X}$, and thus we may assume from start
that $\cF=\cO_{\!X}$. Then, since $X\!\setminus\!U$ is
constructible, corollary \ref{cor_local-depth} applies
and yields the assertion.
\end{proof}

\begin{corollary}\label{cor_depth-cons}
Let $X$ be a locally coherent scheme, $f:X\to S$ a flat,
locally finitely presented morphism, $j:U\to X$ a quasi-compact
open immersion, $\cF$ a reflexive $\cO_{\!X}$-module. Suppose that
\begin{enumerate}
\alphaenu
\item
$\depth_f(x)\geq 1$ for every point $x\in X\!\setminus\! U$,
and
\item
$\depth_f(x)\geq 2$ for every maximal point $\eta$ of $S$
and every $x\in (X\!\setminus\! U)\cap f^{-1}(\eta)$.
\item
$\delta'(s,\cO_{\!S})>0$ for every non-maximal point $s\in S$.
\romanenu
\end{enumerate}
Then the natural morphism\  $\cF\to j_*j^*\cF$\ 
is an isomorphism.
\end{corollary}
\begin{proof}  Since $f$ is flat, corollary
\ref{cor_depth-flat-basechange} and our assumptions
imply that $\delta'(x,\cO_{\!X})\geq 2$ for every
$x\in X\!\setminus\! U$, so the assertion follows
from proposition \ref{prop_extend-rflx}(ii).
\end{proof}

\sset\subsubsection{}\label{subsec_rk-trivias}
Let $X$ be any scheme. Recall that the {\em rank\/} of an
$\cO_{\!X}$-module $\cF$ of finite type, is the upper
semicontinuous function:
$$
\rk\,\cF:X\to\N \qquad
x\mapsto\dim_{\kappa(x)}\cF_x\otimes_{\cO_{\!X,x}}\kappa(x).
$$
Clearly, if $\cF$ is a is locally free $\cO_{\!X}$-module of
finite type, $\rk\,\cF$ is a continuous function on $X$.
The converse holds, provided $X$ is a reduced scheme.
Moreover, if $\cF$ is of finite presentation, $\rk\,\cF$
is a constructible  function and there exists a dense open
subset $U\subset X$ such that $\rk\,\cF$ restricts to a
continuous function on $U$. We denote by $\bPic\,X$ the
full subcategory of $\cO_{\!X}\Mod_\mathrm{lfft}$
consisting of all the objects whose rank is constant
equal to one ({\em i.e.} the {\em invertible\/} $\cO_{\!X}$-modules).
In case $X$ is locally coherent, we shall also consider the category
$\bDiv\,X$ of {\em generically invertible\/} $\cO_{\!X}$-modules,
defined as the full subcategory of $\cO_{\!X}\bRflx$ consisting
of all objects which are locally free of rank one on a dense
open subset of $X$. If $X$ is locally coherent, $\bPic\,X$ is
a full subcategory of $\bDiv X$.

\begin{remark}\label{rem_classy-rflx}
(i)\ \
Let $A$ be any integral domain, and set $X=\Spec\,A$.
Classically, one has a notion of reflexive fractional
ideal of $A$ (see \cite[p.80]{Mat} or example
\ref{ex_integral-doms}). Suppose now that $A$ is also
coherent, in which case we have the notion of reflexive
$A$-module of definition \ref{def_reflex}(iv). We claim
that these two notions overlap on the subclass of
reflexive fractional ideals of finite type: more
precisely, let $\bDiv(A)$ be the full subcategory of
$A\Mod$ whose objects are the finitely generated
reflexive fractional ideals of $A$. Then the essential
image of the functor
\set\begin{equation}\label{eq_two-diff-rflx}
\bDiv(A)\to\cO_{\!X}\Mod
\qquad
M\mapsto M^\sim
\end{equation}
is the category $\bDiv\,X$, and \eqref{eq_two-diff-rflx}
yields an equivalence of $\bDiv(A)$ with the latter category.
Indeed, let $\cF$ be any generically invertible $\cO_{\!X}$-module,
and set $I:=\cF(X)$; if $K$ denotes the field of fractions
of $A$, then $\dim_KI\otimes_AK=1$. Let us then fix a
$K$-linear isomorphism $I\otimes_AK\isom K$, and notice
that the induced $A$-linear map $I\to K$ is injective,
since $\cF$ is $S_1$ (remark \ref{rem_co-repres}). We
may then view $I$ as a finitely generated $A$-submodule of
$K$, and then it is clear that $I$ is a fractional ideal
of $A$. Moreover, on the one hand $\cF^\vee$ is the
coherent $\cO_{\!X}$-module $\Hom_A(I,A)^\sim$; on the
other hand, the natural map 
$\Hom_A(I,A)\to\Hom_K(I\otimes_AK,K)=K$ is injective,
and its image is the fractional ideal $I^{-1}$ (see
\eqref{subsec_fract-ideals}). We easily deduce that $I$ is
a reflexive fractional ideal, and conversely it is easily
seen that every $\cO_{\!X}$-module in the essential image
of \eqref{eq_two-diff-rflx} is generically invertible;
since this functor is also obviously fully faithful, the
assertion follows.

(ii)\ \
Let $A$ be a coherent integral domain, and denote by
$\cohDiv(A)$ the set of all coherent reflexive fractional
ideals of $A$. It is easily seen that $\cohDiv(A)$ is
a submonoid of $\Div(A)$, for the natural monoid structure
introduced in example \ref{ex_integral-doms}. More generally,
if $X$ is a coherent integral scheme, we may define a sheaf
of monoids $\cDiv_X$ on $X$, as follows. First, to any
affine open subset $U\subset X$, we assign the monoid
$\cDiv_X(U):=\cohDiv(\cO_{\!X}(U))$. For each inclusion
$j:U'\subset U$ of affine open subset, notice that the
restriction map $\cO_{\!X}(U)\to\cO_{\!X}(U')$ is a
flat ring homomorphism, hence it gives a flat morphism
of monoids
$\cO_{\!X}(U)\setminus\{0\}\to\cO_{\!X}(U')\setminus\{0\}$.
We then have an induced map of monoids
$\cDiv_X(j):\cDiv_X(U)\to\cDiv_X(U')$, by virtue of
lemma \ref{lem_rflx-rflx}(iv). It is easily seen that the
resulting presheaf $\cDiv_X$ is a sheaf on the site of all
affine open subsets of $X$. By \cite[Ch.0, \S3.2.2]{EGAI},
the latter extends uniquely to a sheaf of monoids on $X$,
which we denote again $\cDiv_X$. We set
$$
\Div(X):=\cDiv_X(X).
$$
If $X$ is normal and locally noetherian (or more generally, if
$X$ is a {\em Krull scheme}, {\em i.e.} $\cO_{\!X}(U)$ is a
Krull ring, for every affine open subset $U\subset X$) then
$\cDiv_X$ is an abelian sheaf, and $\Div(X)$ is an abelian
group (proposition \ref{prop_completely-sat}(i,iii)).
\end{remark}

\sset\subsubsection{}\label{subsec_recall-det}
For future reference, it is useful to recall some
preliminaries concerning the determinant functors defined
in \cite{Kn-Mu}. Let $X$ be a scheme. We denote by
$\grPic\,X$ the category of graded invertible
$\cO_{\!X}$-modules. An object of $\grPic\,X$ is a
pair $(L,\alpha)$, where $L$ is an invertible
$\cO_{\!X}$-module and $\alpha:X\to\Z$ is a continuous
function. A homomorphism $h:(L,\alpha)\to(M,\beta)$
is a homomorphism of $\cO_{\!X}$-modules $h:L\to M$ such that
$h_x=0$ for every $x\in X$ with $\alpha(x)\neq\beta(x)$.
We denote by $\grPic^*X$ the subcategory of
$\grPic\,X$ with the same objects, and whose
morphisms are the isomorphisms in $\grPic\,X$.
Notice that $\grPic\,X$ is a tensor category : the tensor
product of two objects $(L,\alpha)$ and $(M,\beta)$ is
the pair
$(L,\alpha)\otimes(M,\beta):=(L\otimes_{\cO_{\!X}}M,\alpha+\beta)$.
We denote by $\cO_{\!X}\Mod^*_\mathrm{lfft}$ the category whose
objects are the locally free $\cO_{\!X}$-modules of finite
type, and whose morphisms are the $\cO_{\!X}$-linear isomorphisms.
The {\em determinant\/} is the functor:
$$
\det:\cO_{\!X}\Mod^*_\mathrm{lfft}\to\grPic^*X\qquad
F\mapsto(\Lambda^{\rk\,F}_{\cO_{\!X}}F,\rk\,F).
$$
Let $\sD(\cO_{\!X}\Mod)_\mathrm{perf}$ be the category of
perfect complexes of $\cO_{\!X}$-modules; recall that, by
definition, every perfect complex is locally isomorphic
to a bounded complex of locally free $\cO_{\!X}$-modules
of finite type. The category $\sD(\cO_{\!X}\Mod)_\mathrm{perf}^*$
is the subcategory of $\sD(\cO_{\!X}\Mod)_\mathrm{perf}$
with the same objects, and whose morphisms are the
isomorphisms ({\em i.e.} the quasi-isomorphisms of
complexes). The main theorem of chapter 1 of \cite{Kn-Mu}
can be stated as follows.

\begin{lemma}\label{lem_det-fctr} (\cite[Th.1]{Kn-Mu})
With the notation of \eqref{subsec_recall-det} there exists,
for every scheme $X$, an extension of the determinant
functor to a functor:
$$
\sdet:\sD(\cO_{\!X}\Mod)_\mathrm{perf}^*\to\grPic\,X.
$$
These determinant functors commute with every base change.
\qed\end{lemma}

\begin{proposition}\label{prop_Pic-is-Div}
Let $X$ be a regular scheme. Then every reflexive generically
invertible $\cO_{\!X}$-module is invertible.
\end{proposition}
\begin{proof} The question is local on $X$, hence we may
assume that $X$ is affine. Let $\cF$ be a generically
invertible $\cO_{\!X}$-module, and $U\subset X$ a dense
open subset such that $\cF_{|U}$ is invertible. Denote
by $Z_1,\dots,Z_t$ the irreducible components of
$Z:=X\!\setminus\!U$ whose codimension in $X$ equals one,
and for every $i=1,\dots,t$, let $\eta_i$ be the maximal
point of $Z_i$, and set $A_i:=\cO_{\!X,\eta_i}$. Since
$\cF$ is $S_1$ (remark \ref{rem_co-repres}), the stalk
$\cF_{\eta_i}$ is a torsion-free $A_i$-module of finite
type for every $i\leq t$; however, $A_i$ is a discrete
valuation ring, hence $\cF_{\eta_i}$ is a free
$A_i$-module, necessarily of rank one, for $i=1,\dots,t$.
Since $\cF$ is coherent, it follows that there exists
an open neighborhood $U_i$ of $\eta_i$ in $X$, such that
$\cF_{|U_i}$ is a free $\cO_{|U_i}$-module. Hence, we
may replace $U$ by $U\cup U_1\cup\cdots\cup U_t$ and
assume that every irreducible component of $Z$ has
codimension $>1$, therefore $\delta'(x,\cO_{\!X})\geq 2$
for every $x\in Z$. Now, $\cF[0]$ is a perfect complex by
Serre's theorem (\cite[Th.4.4.16]{We}), so the invertible
$\cO_{\!X}$-module $\sdet\,\cF$ is well defined (lemma
\ref{lem_det-fctr}). Let $j:U\to X$ be the open immersion;
in view of proposition \ref{prop_extend-rflx}, we deduce
natural isomorphisms
$$
\cF\isom j_*j^*\cF\isom j_*\sdet(j^*\cF[0])\isom
j_*j^*\sdet\,\cF[0]\stackrel{\sim}{\leftarrow}\sdet\,\cF[0]
$$
and the assertion follows.
\end{proof}

We wish now to introduce a notion of duality better suited
to derived categories of $\cO_{\!X}$-modules (over a scheme $X$).
Hereafter we only carry out a preliminary investigation
of such derived duality -- the full development of which, will
be the task of section \ref{sec_loc-duality}.

\begin{definition}\label{def_dual-complex}
Let $X$ be a locally coherent scheme. A complex $\omega^\bullet$
in $\sD^b(\cO_{\!X}\Mod)_\coh$ (notation of
\eqref{sec_various-O-mod}) is called {\em dualizing\/}
if it fulfills the following two conditions :
\begin{enumerate}
\alphaenu
\item
The functor :
$$
\cD:\sD(\cO_{\!X}\Mod)^o\to\sD(\cO_{\!X}\Mod)\qquad
C^\bullet\mapsto R\cHom^\bullet_{\cO_{\!X}}(C^\bullet,\omega^\bullet)
$$
restricts to a {\em duality functor} :
$\cD:\sD^b(\cO_{\!X}\Mod)^o_\coh\to\sD^b(\cO_{\!X}\Mod)_\coh$.
\item
The natural transformation :
$\eta_{C^\bullet}:C^\bullet\to\cD\circ\cD(C^\bullet)$
restricts to a {\em biduality isomorphism\/} of functors
on the category $\sD^b(\cO_{\!X}\Mod)_\coh$.
\romanenu
\end{enumerate}
\end{definition}

\begin{remark}\label{rem_dual-complex}
Let $X$ be a locally coherent scheme, $\omega^\bullet$ an object
of $\sD^b(\cO_{\!}\Mod)_\coh$, and define the functor $\cD$ as in
definition \ref{def_dual-complex}.

(i)\ \
A standard {\em d\'evissage\/} argument shows that
$\omega^\bullet$ is dualizing on $X$ if and only if
$\cD(\cF[0])$ lies in $\sD^b(\cO_{\!X}\Mod)_\coh$, and
the biduality map $\cF[0]\to\cD\circ\cD(\cF[0])$ is an
isomorphism for every coherent $\cO_{\!X}$-module $\cF$.

(ii)\ \
Suppose that the complex $\omega^\bullet$ is dualizing
on $X$. From the natural identification
$\omega^\bullet\isom R\cHom^\bullet_{\cO_{\!X}}(\cO_{\!X}[0],\omega^\bullet)$,
we deduce a natural isomorphism
$$
R\cHom^\bullet_{\cO_{\!X}}
(\omega^\bullet,\omega^\bullet)\isom\cO_{\!X}[0]
\qquad
\text{in $\sD^b(\cO_{\!X}\Mod)$}.
$$
\end{remark}

\begin{example}\label{ex_reg-is-Gorenstein}
Suppose that $X$ is a noetherian regular scheme
({\em i.e.} all the stalks $\cO_{\!X,x}$ are regular rings). In
light of Serre's theorem \cite[Th.4.4.16]{We}, every object of
$\sD^b(\cO_{\!X}\Mod)_\coh$ is a perfect complex. It follows
easily that the complex $\cO_{\!X}[0]$ is dualizing. For more
general schemes, the structure sheaf does not necessarily work,
and the existence of a dualizing complex is a delicate issue.
On the other hand, one may ask to what extent a complex is
determined by the properties (a) and (b) of definition
\ref{def_dual-complex}. Clearly, if $\omega^\bullet$ is
dualizing on $X$, then so is any other complex of the form
$\omega^\bullet\otimes_{\cO_{\!X}}\cL$, where $\cL$ is an
invertible $\cO_{\!X}$-module. Also, any shift of $\omega^\bullet$
is again dualizing. Conversely, the following proposition
\ref{prop_unique-dual} says that any two dualizing complexes
are related in such manner, up to quasi-isomorphism.
\end{example}

\begin{lemma}\label{lem_two-duals}
Let $(X,\cO_{\!X})$ be any locally ringed space, and
$P^\bullet$ and $Q^\bullet$ two objects of
$\sD^-(\cO_{\!X}\Mod)$ with a quasi-isomorphism :
$$
P^\bullet\derotimes_{\cO_{\!X}}Q^\bullet\isom\cO_{\!X}[0].
$$
Then there exist an invertible $\cO_{\!X}$-module $\cL$,
a continuous function $\sigma:|X|\to\Z$ and quasi-isomorphisms:
$$
P^\bullet\isom\cL[\sigma] \qquad\text{and}\qquad
Q^\bullet\isom\cL^{-1}[-\sigma].
$$
\end{lemma}
\begin{proof} It suffices to verify that, locally on $X$,
the complexes $P^\bullet$ and $Q^\bullet$ are of the required
form; indeed in this case $X$ will be a disjoint union of open
sets $U_n$ on which $H^\bullet P^\bullet$ is concentrated in
degree $n$ and $H^\bullet Q^\bullet$ is concentrated in degree $-n$.
Note that $\cO_{\!X}\Mod$ is equivalent to the product of the
categories $\cO_{\!U_n}\Mod$ and that the derived category of
a product of abelian categories is the product of  the derived
categories of the factors.

\begin{claim}\label{cl_OK-local} Let $A$ be a (commutative) local
ring, $K^\bullet$ and $L^\bullet$ two objects of $\sD^-(A\Mod)$ with
a quasi-isomorphism
\set\begin{equation}\label{eq_give-it}
K^\bullet\derotimes_AL^\bullet\isom A[0].
\end{equation}
Then there exist $s\in\Z$ and quasi-isomorphisms :
$$
K^\bullet\isom A[s]
\qquad
L^\bullet\isom A[-s].
$$
\end{claim}
\begin{pfclaim} Set :
$$
i_0:=\max(i\in\Z~|~(H^iK^\bullet)\neq 0)
\quad\text{and}\quad
j_0:=\max(i\in\Z~|~(H^iL^\bullet)\neq 0).
$$
We may assume that $K^i=0$ for every $i>i_0$, and that $L^\bullet$
is a bounded above complex of free $A$-modules. Then we may find
a filtered system $(K^\bullet_\lambda~|~\lambda\in\Lambda)$ of
complexes of $A$-modules bounded from above, such that 
\begin{itemize}
\item
$H^{i_0}(K^\bullet_\lambda)$ is a finitely generated $A$-module, for
every $\lambda\in\Lambda$;
\item
the colimit of the system $(K^\bullet_\lambda~|~\lambda\in\Lambda)$
(in the category of complexes of $A$-modules) is isomorphic to
$K^\bullet$.
\end{itemize}
From \eqref{eq_give-it} we get an isomorphism
$H^0(K^\bullet\otimes_AL^\bullet)\isom A$, and it follows easily that
the natural map
$H^0(K^\bullet_\lambda\otimes_AL^\bullet)\to H^0(K^\bullet\otimes_AL^\bullet)$
is surjective for some $\lambda\in\Lambda$. For such $\lambda$, we
may then find a morphism in $\sD^-(A\Mod)$ :
\set\begin{equation}\label{eq_give-it-again}
A[0]\to K^\bullet_\lambda\derotimes_AL^\bullet
\end{equation}
whose composition with the natural map
$K^\bullet_\lambda\derotimes_AL^\bullet\to K^\bullet\derotimes_AL^\bullet$
is the inverse of \eqref{eq_give-it}. Hence
$$
\eqref{eq_give-it-again}\derotimes_AK^\bullet:
K^\bullet\to(K^\bullet_\lambda\derotimes_AL^\bullet)\derotimes_AK^\bullet
\isom
K^\bullet_\lambda\derotimes_A(L^\bullet\derotimes_AK^\bullet)\isom
K^\bullet_\lambda
$$
is a right inverse of the natural morphism $K^\bullet_\lambda\to K^\bullet$
in $\sD^-(A\Mod)$. Especially, the induced map
$H^{i_0}K^\bullet_\lambda\to H^{i_0}K^\bullet$ is surjective, {\em i.e.}
$H^{i_0}K^\bullet$ is a finitely generated $A$-module. Likewise, we
see that $H^{j_0}L^\bullet$ is a finitely generated $A$-module.
Now, notice that 
$$
H^k(K^\bullet\derotimes_AL^\bullet)\simeq\left\{
\begin{array}{ll}
H^{i_0}K^\bullet\otimes_AH^{j_0}L^\bullet
& \text{for $k=i_0+j_0$} \\
0 & \text{for $k>i_0+j_0$}.
\end{array}\right.
$$
From Nakayama's lemma it follows easily that
$H^{i_0+j_0}(K^\bullet\derotimes_AL^\bullet)\neq 0$,
and then our assumptions imply that $i_0+j_0=0$ and
$H^{i_0}K^\bullet\otimes_AH^{j_0}L^\bullet\simeq A$.
One deduces easily that
$H^{i_0}K^\bullet\simeq A\simeq H^{j_0}L^\bullet$
(see {\em e.g.} \cite[Lemma 4.1.5]{Ga-Ra}). Furthermore, we
can find a complex $K_{\!1}^\bullet$ in $\sD^{<i_0}(A\Mod)$
(resp. $L_1^\bullet$ in $\sD^{<j_0}(A\Mod)$) such that :
$$
K^\bullet\simeq A[-i_0]\oplus K^\bullet_{\!1}
\qquad\text{(resp. $L^\bullet\simeq A[-j_0]\oplus L^\bullet_1$)}
$$
whence a quasi-isomorphism :
$$
\phi:A[0]\isom K^\bullet\derotimes_AL^\bullet\isom
A[0]\oplus K^\bullet_{\!1}[-j_0]\oplus L^\bullet_1[-i_0]
\oplus(K^\bullet_{\!1}\derotimes_AL^\bullet_1).
$$
However, by construction $\phi^{-1}$ restricts to an isomorphism on
the direct summand $A[0]$, therefore
$K^\bullet_{\!1}\simeq 0\simeq L^\bullet_1$ in $\sD(A\Mod)$,
and the claim follows.
\end{pfclaim}

Now, for any point $x\in X$, let $i_x:\{x\}\to X$ be the inclusion
map, and set $K^\bullet_x:=i^*_xK^\bullet$ for every complex
$K^\bullet$ of $\cO_{\!X}$-modules. Notice that if $K^\bullet$
is a complex of flat $\cO_{\!X}$-modules, then $K^\bullet_x$
is a complex of flat $\cO_{\!X,x}$-modules
(\cite[Exp.V, Prop.1.6(1)]{SGA4-2}), therefore the rule
$K^\bullet\mapsto K^\bullet_x$ yields a well-defined functor
$\sD^-(\cO_{\!X}\Mod)\to\sD^-(\cO_{\!X,x}\Mod)$,
and moreover we have a natural isomorphism
\set\begin{equation}\label{eq_loc-derotimes}
K^\bullet_x\derotimes_{\cO_{\!X,x}}L^\bullet_x\isom
(K^\bullet\derotimes_{\cO_{\!X}}L^\bullet)_x
\end{equation}
for every objects $K^\bullet,L^\bullet$ of $\sD^-(\cO_{\!X}\Mod)$.
Especially, under the current assumptions, and in view of claim
\ref{cl_OK-local}, we may find $s\in\Z$, and an isomorphism of
$\cO_{\!X,x}$-modules
\set\begin{equation}\label{eq_mapsto-one}
H^sP^\bullet_x\otimes_{\cO_{\!X,x}}H^{-s}Q^\bullet_x\isom
H^0(P^\bullet_x\derotimes_{\cO_{\!X,x}}Q^\bullet_x)\isom\cO_{\!X,x}.
\end{equation}
Since $\cO_{\!X,x}$ is local, we may thus find $a_x\in H^sP^\bullet_x$
and $b_x\in H^{-s}Q^\bullet_x$ such that \eqref{eq_mapsto-one} maps
$a_x\otimes b_x$ to $1$. Then $a_x$ and $b_x$ extend to local sections
$$
a\in\Gamma(U,\Ker(d:P^s\to P^{s-1}))
\qquad
b\in\Gamma(U,\Ker(d:Q^{-s}\to P^{-s-1}))
$$
on some neighborhood $U\subset X$ of $x$, and after shrinking
$U$, we may assume that $a\otimes b$ gets mapped to $1$,
under the induced morphism of $\cO_{\!U}$-modules
$$
H^sP^\bullet_{|U}\otimes_{\cO_{\!U}}H^{-s}Q^\bullet_{|U}\to
H^0(P^\bullet\derotimes_{\cO_{\!X}}Q^\bullet)_{|U}\isom\cO_{\!U}.
$$
Then we obtain a well defined morphism in $\sD^-(\cO_{\!U}\Mod)$
$$
\phi:\cO_{\!U}[s]\to P^\bullet_{|U}
\qquad
\text{(resp.\quad $\psi:\cO_{\!U}[-s]\to Q^\bullet_{|U}$)}
$$
by the rule : $t\mapsto t\cdot a$ (resp. $t\mapsto t\cdot b$)
for every local section $t$ of $\cO_{\!U}$. Again by claim
\ref{cl_OK-local} and \eqref{eq_loc-derotimes} we deduce that
$\phi_y:\cO_{\!U,y}[s]\to P^\bullet_y$ is a quasi-isomorphism
for every $y\in U$ (and likewise for $\psi_y$); {\em. i.e.}
$\phi$ and $\psi$ are the sought isomorphisms in
$\sD^-(\cO_{\!U}\Mod)$.
\end{proof}

\begin{proposition}\label{prop_unique-dual}
Suppose that $\omega_1^\bullet$ and $\omega_2^\bullet$ are
two dualizing complexes for the locally coherent scheme $X$.
Then there exist an invertible $\cO_{\!X}$-module $\cL$ and
a continuous function $\sigma:|X|\to\Z$ such that
$$
\omega_2^\bullet\simeq\omega_1^\bullet\otimes_{\cO_{\!X}}\cL[\sigma]
\qquad\text{in $\sD^b(\cO_{\!X}\Mod)_\coh$}.
$$
\end{proposition}
\begin{proof} Denote by $\cD_1$ and $\cD_2$ the duality functors
associated with $\omega_1$ and  respectively $\omega_2$.
By assumption, we can find complexes $P^\bullet,Q^\bullet$
in $\sD^b(\cO_{\!X}\Mod)_\coh$ such that
$\omega_2\simeq\cD_1(P^\bullet)$ and $\omega_1\simeq\cD_2(Q^\bullet)$,
and therefore
$$
\cD_2(\cF^\bullet)\simeq
R\cHom^\bullet_{\cO_{\!X}}(\cF^\bullet,\cD_1(P^\bullet))\simeq
\cD_1(\cF^\bullet\derotimes_{\cO_{\!X}}P^\bullet)
$$
for every object $\cF^\bullet$ of $\sD^b(\cO_{\!X}\Mod)_\coh$
(\cite[Th.10.8.7]{We}).

\begin{claim}\label{cl_der-bddness} Let $C^\bullet$ be an object of
$\sD^-(\cO_{\!X}\Mod)_\coh$, such that $\cD_1(C^\bullet)$ is in
$\sD^b(\cO_{\!X}\Mod)$. Then $C^\bullet$ is in $\sD^b(\cO_{\!X}\Mod)_\coh$.
\end{claim}
\begin{pfclaim} For given $m,n\in\N$, the natural
maps: $\tau^{\leq-n}C^\bullet\xrightarrow{\alpha}
C^\bullet\xrightarrow{\beta}\tau^{\geq-m}C^\bullet$
induce morphisms
$$
\cD_1(\tau^{\geq-m}C^\bullet)\xrightarrow{\cD_1(\beta)}
\cD_1(C^\bullet) \xrightarrow{\cD_1(\alpha)}\cD_1(\tau^{\leq-n}C^\bullet)
$$
Say that $\omega^\bullet\simeq\tau_{\geq a}\omega^\bullet$ for
some integer $a\in\N$. Then $\cD_1(\tau^{\leq-n}C^\bullet)$
lies in $\sD^{\geq n+a}(\cO_{\!X}\Mod)$. Since $\cD_1(C^\bullet)$
is bounded, it follows that $\cD_1(\beta)=0$ for $n$ large enough.
Consider now the commutative diagram :
$$
\xymatrix{
\tau^{\geq-n}C^\bullet \ar[rr]^{\beta\circ\alpha} \ar[d] & &
\tau^{\geq-m}C^\bullet \ar[d]^\eta \\
\cD_1\circ\cD_1(\tau^{\geq-n}C^\bullet)
\ar[rr]^{\cD_1\circ\cD_1(\beta\circ\alpha)} & &
\cD_1\circ\cD_1(\tau^{\geq-m}C^\bullet)
}$$
Since $\tau^{\geq-m}C^\bullet$ is a bounded complex, $\eta$
is an isomorphism in $\sD(\cO_{\!X}\Mod)$, so $\beta\circ\alpha=0$
whenever $n$ is large enough. Clearly this means that
$C^\bullet$ is bounded, as claimed.
\end{pfclaim}

Applying claim \ref{cl_der-bddness} to
$C^\bullet:=\cF^\bullet\derotimes_{\cO_X}P^\bullet$
(which is in $\sD^-(\cO_{\!X}\Mod)_\coh$, since $X$ is coherent)
we see that the latter is a bounded complex, and by
reversing the roles of $\omega_1$ and $\omega_2$ it follows
that the same holds for $\cF^\bullet\derotimes_{\cO_{\!X}}Q^\bullet$.
We then deduce isomorphisms in $\sD^b(\cO_{\!X}\Mod)_\coh$
$$
\cD_1\circ\cD_2(\cF^\bullet)\simeq
\cF^\bullet\derotimes_{\cO_{\!X}}P^\bullet
\qquad\text{and}\qquad
\cD_2\circ\cD_1(\cF^\bullet)\simeq
\cF^\bullet\derotimes_{\cO_{\!X}}Q^\bullet.
$$
Letting $\cF^\bullet:=\cO_{\!X}[0]$ we derive :
$$
\cO_{\!X}[0]\simeq\cD_2\circ\cD_1\circ\cD_1\circ\cD_2(\cO_{\!X}[0])\simeq
\cD_2\circ\cD_1(P^\bullet)\simeq P^\bullet\derotimes_{\cO_{\!X}}Q^\bullet.
$$
Then lemma \ref{lem_two-duals} says that $P^\bullet\simeq\cE[\tau]$
for an invertible $\cO_{\!X}$-module $\cE$ and a continuous function
$\tau:|X|\to\Z$. Consequently :
$\omega_2\simeq\cE^\vee[-\tau]\otimes_{\cO_{\!X}}\omega_1$, so the
proposition holds with $\cL:=\cE^\vee$ and $\sigma:=-\tau$.
\end{proof}

\begin{lemma}\label{lem_transit-dual}
Let $f:X\to Y$ be a morphism of locally coherent schemes,
and $\omega_Y^\bullet$ a dualizing complex on $Y$. We have :
\begin{enumerate}
\item
If $f$ is finite and finitely presented, then
$f^!\omega_Y^\bullet$ is dualizing on $X$.
\item
If\/ $Y$ coherent and $f$ is an open immersion, then
$f^*\omega^\bullet_Y$ is dualizing on $X$.
\end{enumerate}
\end{lemma}
\begin{proof}(i): Denote by
$\bar f:(X,\cO_{\!X})\to(Y,f_*\cO_{\!X})$ the morphism of
ringed spaces deduced from $f$. For any object $C^\bullet$
of $\sD^-(\cO_{\!X}\Mod)_\coh$ we have natural isomorphisms :
$$
\begin{aligned}
\cD(C^\bullet):=R\cHom_{\cO_X}(C^\bullet,f^!\omega_Y^\bullet) &
\isom R\cHom_{\cO_{\!X}}(C^\bullet,
\bar f{}^*R\cHom_{\cO_Y}(f_*\cO_{\!X},\omega_Y^\bullet)) \\
& \isom\bar f{}^*R\cHom_{f_*\cO_{\!X}}
(f_*C^\bullet,R\cHom_{\cO_Y}(f_*\cO_{\!X},\omega_Y^\bullet)) \\
& \isom\bar f{}^*R\cHom_{\cO_Y}(f_*C^\bullet,\omega_Y^\bullet).
\end{aligned}
$$
Hence, if $C^\bullet$ is in $\sD^b(\cO_{\!X}\Mod)_\coh$, the
same holds for $\cD(C^\bullet)$; especially, this applies to
$f^!\omega_Y\simeq\cD(\cO_{\!X}[0])$, and we can compute :
$$
\begin{aligned}
\cD\circ\cD(C^\bullet) & \isom\bar f{}^*R\cHom_{\cO_Y}
(R\cHom_{\cO_Y}(f_*C^\bullet,\omega_Y^\bullet),\omega_Y^\bullet) \\
& \isom\bar f{}^*f_*C^\bullet \isom C^\bullet
\end{aligned}
$$
and by inspecting the definitions, one verifies that the
resulting natural transformation $C^\bullet\to\cD\circ\cD(C^\bullet)$
is the biduality isomorphism. The claim follows.

(ii): Since $\cO_Y$ is coherent, the natural map
$$
f^*R\cHom_{\cO_{Y}}(\cG[0],\omega^\bullet_Y)\to
R\cHom_{\cO_{\!X}}(f^*\cG[0],f^*\omega^\bullet_Y)
$$
is an isomorphism in $\sD^b(\cO_{\!X}\Mod)_\coh$, for
every $\cO_Y$-module $\cG$. Then the assertion follows
easily from lemma \ref{lem_extend-cohs}(ii) and remark
\ref{rem_dual-complex}(i).
\end{proof}

\begin{proposition}\label{prop_Stanley-criterion}
Let $A$ be a noetherian local ring, $\kappa$ its residue
field, $M^\bullet$ a bounded below complex of $A$-modules of
finite type, and set $X:=\Spec\,A$. We have :
\begin{enumerate}
\item
The following conditions are equivalent 
\begin{enumerate}
\item
There exists $c\in\Z$ such that
$R\Hom^\bullet_A(\kappa[0],M^\bullet)\simeq\kappa[c]$.
\item
The complex of $\cO_{\!X}$-modules $M^{\bullet\sim}$
arising from $M^\bullet$ is dualizing on $X$.
\end{enumerate}
\item
If the equivalent conditions of\/ {\em(i)} hold, then
$M^\bullet$ has finite injective dimension. 
\end{enumerate}
\end{proposition}
\begin{proof}(i): We show first that (a)$\Rightarrow$(b).
To this aim, in view of remark \ref{rem_dual-complex}(i) and
corollary \ref{cor_Ext-loc=glob}(ii), it suffices to check :

\begin{claim}\label{cl_new-name}
Assume (a). Then, for every finitely generated $A$-module $N$,
we have :
\begin{enumerate}
\item
$D^\bullet(N):=R\Hom^\bullet_A(N[0],M^\bullet)\in\sD^{\leq-c}(A\Mod)$
(notation of \eqref{sec_brutal-truncate}).
\item
The natural map
$$
N[0]\to DD^\bullet(N):=R\Hom^\bullet_A(D^\bullet(N),M^\bullet)
$$
is an isomorphism.
\end{enumerate}
\end{claim}
\begin{pfclaim} We first show the claim for $N=\kappa$,
in which case (i) holds by assumption. To check (ii), let
$M^\bullet\isom I^\bullet$ be a resolution consisting of a
bounded below complex of injective $A$-modules, so that
$$
D^\bullet(\kappa[0],M^\bullet)\isom\Hom^\bullet_A(\kappa[0],I^\bullet)
\isom I^\bullet[\fm]
$$
where $\fm\subset A$ is the maximal ideal, and $I^k[\fm]$ denotes
the submodule of $\fm$-torsion elements in $I^k$, for every $k\in\Z$.
Under these isomorphisms, it is easily seen that the biduality map
of (ii) is identified with the unique one
$$
\kappa\to H:=\rHot_{A\Mod}(I^\bullet[\fm],I^\bullet)
$$
that sends $1\in\kappa$ to the inclusion map
$j^\bullet:I^\bullet[\fm]\to I^\bullet$ (details left to the
reader). Now, it is clear that any morphism
$I^\bullet[\fm]\to I^\bullet$ in $\Hot(A\Mod)$ factors through
$j^\bullet$, and on the other hand, our assumption implies that
$H\simeq\kappa$. Hence, pick a morphism
$f^\bullet:I^\bullet[\fm]\to I^\bullet$ representing a generator
for the $A$-module $H$, and write $f^\bullet=j^\bullet\circ g^\bullet$
for some endomorphism $g^\bullet$ of $I^\bullet[\fm]$; in other
words, $f^\bullet$ is the image of $j^\bullet$ under the
$A$-linear map
$$
\rHot_{A\Mod}(g^\bullet,I^\bullet):H\to H
$$
so the class of $j^\bullet$ cannot vanish in $H$, and (ii)
follows in this case.

Next, we shall argue by induction on $d:=\dim\,\Supp\,N$.
If $d=0$, then $N$ is an $A$-module of finite length, in
which case we argue by induction on the length $l$ of $N$.
If $l=1$, we have $N\simeq\kappa$, so the assertions are
already known. Suppose $l>1$, and that both (i) and (ii)
are already known for all $A$-modules of length $<d$; we
may find an $A$-submodule $N'\subset N$ such that both
$N'$ and $N'':=N/N'$ have length $<d$. From the inductive
assumption for $N'$ and $N''$, and the induced distiguished
triangle
$$
D^\bullet(N'')\to D^\bullet(N)\to D^\bullet(N')\to D^\bullet(N'')[1]
$$
we deduce that (i) holds for $N$. Likewise, since (ii) is
known for both $N'$ and $N''$, using the $5$-lemma we deduce
easily that the same holds also for $N$.

Lastly, suppose that $d>0$, and both (i) and (ii) are already
known for all $A$-modules of finite type whose support has
dimension $<d$. Let $N':=\Gamma_{\!\{\fm\}}N$; both (i) and (ii)
are already known for $N'$, so the same {\em d\'evissage}
argument as in the foregoing reduces to showing the claim
for $N/N'$, {\em i.e.} we may assume that $\fm\notin\Ass\,N$.
Thus, let $t\in\fm$ be any element such that the scalar
multiplication map $t\cdot\one_N$ is injective, so we have
a short exact sequence
$$
0\to N\xrightarrow{\ t^k\ } N\to N_k:=N/t^kN\to 0
\qquad
\text{for every $k>0$}.
$$
Notice that $\dim\Supp\,N_k<d$ : indeed, if $\fp$ is a
minimal element of $\Supp\,N$, then $\fp\in\Ass\,N$
(\cite[Th.6.5(iii)]{Mat}), hence $t\notin\fp$, and
therefore $\fp\notin\Supp\,N_k$. By inductive assumption,
both (i) and (ii) hold for $N_k$ (for every $k>0$), and by
considering the induced distinguished triangle
$$
D^\bullet(N_k)\to D^\bullet(N)\xrightarrow{\ t^k\ }D^\bullet(N)
\to D^\bullet(N_k)[1]
$$
we deduce that scalar multiplication by $t$ is an epimorphism
on $H^jD^\bullet(N)$ whenever $j\geq-c$; but the latter is an
$A$-module of finite type, so it must vanish, by Nakayama'
lemma, and we conclude that (i) holds for $N$. Furthermore,
by the same token we get an exact sequence
$$
H^iDD^\bullet(N)\xrightarrow{\ t^k\ }H^iDD^\bullet(N)\to 0
\qquad
\text{for every $i\neq 0$}
$$
whence $H^iDD(N)=0$ for $i\neq 0$, again by Nakayama's lemma.
Lastly, for every $k>0$ consider the ladder with exact rows :
$$
\xymatrix{
0 \ar[r] & N \ar[r]^-{t^k} \ar[d]_\alpha & N \ar[r] \ar[d]^\alpha &
N_k \ar[r] \ar[d] & 0 \\
0 \ar[r] & H^0DD^\bullet(N) \ar[r]^-{t^k} & H^0DD^\bullet(N) \ar[r] &
H^0DD^\bullet(N_k) \ar[r] & 0
}$$
whose right-most vertical arrow is an isomorphism, by inductive
assumption. A simple diagram chase then yields
$$
H^0DD^\bullet(N)=t^k\cdot H^0DD^\bullet(N)+\alpha(N)
$$
so $\alpha$ is surjective, by Nakayama's lemma. To show
the injectivity of $\alpha$, let $x\in N$ be any non-zero
element, and choose $k>0$ such that $x\notin t^kN$, so the
image of $x$ does not vanish in $N_k$, whence necessarily
$\alpha(x)\neq 0$, and the claim follows.
\end{pfclaim}

(b)$\Rightarrow$(a): Let $i:\Spec\,\kappa\to X$ be the
closed immersion; from corollary \ref{cor_Ext-loc=glob}(ii)
we know that $i^!M^{\bullet\sim}$ is the complex arising from
the complex of $\kappa$-modules $R\Hom_A(\kappa[0],M^\bullet)$.
On the other hand, $i^!M^{\bullet\sim}$ is dualizing on
$\Spec\,\kappa$, by virtue of lemma \ref{lem_transit-dual}(i);
obviously, $\kappa[0]^\sim$ is also dualizing on the same
scheme, so the assertion follows from proposition
\ref{prop_unique-dual}.

(ii): We notice more precisely :

\begin{claim}\label{cl_obvious-shift}
Suppose that the conditions of (i) hold for $M^\bullet$, and
moreover that $H^0M^\bullet\neq 0$ and $H^iM^\bullet=0$ for
every $i<0$. Then $c\leq 0$, and the injective dimension of
$M^\bullet$ equals $-c$.
\end{claim}
\begin{pfclaim} Under the assumptions of the claim, it is
clear that $c\leq 0$. Pick a quasi-isomorphism
$M^\bullet\isom I^\bullet$, where $I^\bullet$ is an object
of $\sC^{[0,-c]}(A\Mod)$ such that $I^j$ is an injective
$A$-module for every $j<-c$. Notice that $I^0\neq 0$, since
$H^0M^\bullet$ does not vanish. Let $N$ be any $A$-module;
a standard argument shows that
$$
\Ext^1_A(N,I^c)\simeq
R^{1-c}\Hom_A^\bullet(N[0],M^\bullet)
$$
and the latter vanishes if $N$ is finitely generated, by
virtue of claim \ref{cl_new-name}(i). It follows that $I^c$
is injective as well (see {\em e.g.} \cite[Th.B3]{Mat}),
so the injective dimension of $M^\bullet$ is $\leq -c$.
But it is also $\geq-c$, due to condition (i.a).
\end{pfclaim}

After replacing $M^\bullet$ by $M^\bullet[a]$ for a suitable
$a\in\Z$, the assumptions of claim \ref{cl_obvious-shift}
can obviously be fulfilled, whence (ii).
\end{proof}

Let us recall the following :

\begin{definition}\label{def_Gorenstein}
Let $X$ be a noetherian scheme.

(i)\ \
We say that $X$ is {\em Gorenstein}, if $\cO_{\!X}[0]$
is dualizing on $X$.

(ii)\ \
If $X=\Spec\,A$ is affine and Gorenstein, we also say
that $A$ is a {\em Gorenstein ring}.
\end{definition}

\begin{remark}\label{rem_Gorenstein}
(i)\ \
In view of example \ref{ex_reg-is-Gorenstein}, every
regular noetherian ring is Gorenstein.

(ii)\ \
Let $A$ be any local noetherian ring, and $\kappa$
the residue field of $A$. In light of proposition
\ref{prop_Stanley-criterion}(i) (and of corollary
\ref{cor_Ext-loc=glob}(ii)), we see that $A$ is Gorenstein
if and only if there exists $c\in\N$ such that
$R\Hom_A^\bullet(\kappa[0],A[0])\simeq\kappa[c]$.
\end{remark}

\sset\subsubsection{}\label{subsec_pull-back-dualize}
Let $A\to B$ be a flat homomorphism of local noetherian rings,
$\kappa_A$ (resp. $\kappa_B$) the residue field of $A$ (resp.
of $B$), and $K^\bullet$ a bounded complex of $A$-modules
of finite type. Set $X:=\Spec\,A$, $Y:=\Spec\,B$, denote
by $f:Y\to X$ the resulting morphism of local schemes, and
let $K^{\bullet\sim}$ be the complex of $\cO_{\!X}$-modules
arising from $K^\bullet$. We have :

\begin{proposition}\label{prop_Goren-pullback}
In the situation of \eqref{subsec_pull-back-dualize}, the
following conditions are equivalent :
\begin{enumerate}
\alphaenu
\item
$f^*K^{\bullet\sim}$ is dualizing on $Y$.
\item
$K^{\bullet\sim}$ is dualizing on $X$ and $B\otimes_A\kappa_A$
is a Gorenstein ring.
\end{enumerate}
\end{proposition}
\begin{proof} Set $C:=B\otimes_A\kappa_A$, $F:=\Spec\,C$ and
$Z:=\Spec\,\kappa_B$, denote $i_1:Z\to F$ and $i_2:F\to Y$
the closed immersions, and let $i:=i_2\circ i_1$.
In light of proposition \ref{prop_Stanley-criterion}(i)
(and of corollary \ref{cor_Ext-loc=glob}(ii)) we see that
(a) holds if and only if there exists $c\in\Z$ such that
$i^!(f^*K^{\bullet\sim})\simeq\cO_{\!Z}[c]$. However,
recall that $i^!\simeq i_1^!\circ i^!_2$ (proposition
\ref{prop_sharp-flat}(i)); invoking again proposition
\ref{prop_Stanley-criterion}(i) we deduce that (a) holds
if and only if $i_2^!(f^*K^{\bullet\sim})$ is dualizing
on $F$. The latter is the complex of $\cO_{\!F}$-modules
arising from
\set\begin{equation}\label{eq_woppity-doo}
R\Hom_B^\bullet(C[0],B\otimes_AK^\bullet)\simeq B\otimes_AL^\bullet
\quad\text{where}\quad
L^\bullet:=R\Hom_A^\bullet(\kappa_A[0],K^\bullet)
\end{equation}
(proposition \ref{prop_replace-prop.12.3.5}(ii) and
corollary \ref{cor_Ext-loc=glob}(ii)). Hence, suppose that
(b) holds; from proposition \ref{prop_Stanley-criterion}(i)
we conclude that $i_2^!(f^*K^{\bullet\sim})\simeq\cO_F[c]$
for some $c\in\Z$; the latter is dualizing on $Z$ by
assumption, whence (a). Conversely, suppose that (a) holds;
from \eqref{eq_woppity-doo} we see that $B\otimes_AL^\bullet$
is a bounded complex of $C$-modules whose cohomology is a
free $C$-module in every degree. We recall the following :

\begin{claim}\label{cl_decompose-in-der}
Let $R$ be any ring, $M^\bullet$ a bounded complex of $R$-modules,
such that $H^iM^\bullet$ is a projective $R$-module, for every
$i\in\Z$. Then we have an isomorphism
$$
M^\bullet\isom\bigoplus_{i\in\Z}(H^iM^\bullet)[-i]
\qquad
\text{in $\sD(R\Mod)$}.
$$
\end{claim}
\begin{pfclaim} We argue by induction on the cardinality $c_M$
of $S_M:=\{i\in\Z~|~H^iM^\bullet\neq 0\}$. If $c_M=0$, the assertion
is obvious. Suppose that $c_M>0$, and that the assertion is known
for every complex $P^\bullet$ as in the claim, with $c_P<c_M$.
Let $a\in\Z$ be the largest element of $S_M$; after replacing
$M^\bullet$ by its truncation $\tau^{\leq a}M^\bullet$, we may
assume that $M^i=0$ for every $i>a$. Then
$H^aM^\bullet\simeq\Coker\,(d^{a-1}:M^{a-1}\to M^a)$, and since
$H^aM^\bullet$ is a projective $R$-module, we may find an
$R$-linear isomorphism $M\isom(H^aM^\bullet)\oplus\Img\,d^{a-1}$.
Denote by $P^\bullet$ the complex such that $P^i:=M^i$ for every
integer $i<a$, and such that $P^a:=\Img\,d^{a-1}$, with differentials
induced by those of $M^\bullet$; then
$M^\bullet\simeq(H^aM^\bullet)[-a]\oplus P^\bullet$ in $\sD(R\Mod)$,
and $c_P<c_M$, so the assertion is known for $P^\bullet$, and the
claim follows.
\end{pfclaim}

In light of claim \ref{cl_decompose-in-der} and remark
\ref{rem_dual-complex}(ii), we easily see that
$B\otimes_AL^\bullet\simeq C[a]$ in $\sD^b(C\Mod)$ for some
$a\in\Z$ (details left to the reader), and then clearly
$L^\bullet\simeq k_A[c]$ in $\sD^b(A\Mod)$. Taking into
account proposition \ref{prop_Stanley-criterion}(i),
assertion (b) follows.
\end{proof}

\begin{proposition}\label{prop_fishy}
Let $X$ be a noetherian scheme, $\cK^\bullet$ an object
of $\sD^b(\cO_{\!X}\Mod)_\coh$. The following conditions
are equivalent :
\begin{enumerate}
\alphaenu
\item
$\cK^\bullet$ is dualizing on $X$.
\item
$\cK^\bullet(x)$ is dualizing on $X(x)$, for every $x\in X$
(notation of definition {\em\ref{def_strict-loc}(iii)}).
\end{enumerate}
\end{proposition}
\begin{proof}(b)$\Rightarrow$(a): We notice the following :

\begin{claim}\label{cl_fin-inj-dim}
Let $X$ be a noetherian scheme, $\cK^\bullet$ an object
of $\sD^b(\cO_{\!X}\Mod)_\coh$, and $x\in X$ any point.
Suppose that the complex of $\cO_{\!X(x)}$-modules
$\cK^\bullet(x)$ has finite injective dimension. Then,
for every $\cF^\bullet\in\sD^b(\cO_{\!X,x}\Mod)_\coh$ there
exists an open neighborhood $U$ of $x$ in $X$ such that
$R\cHom^\bullet_{\cO_{\!X}}(\cF^\bullet,\cK^\bullet)_{|U}$ lies in
$\sD^b(\cO_{\!U}\Mod)_\coh$.
\end{claim}
\begin{pfclaim} The question is local on $X$, so we may
assume that $X$ is affine, say $X=\Spec\,A$ for a noetherian
ring $A$; also, by theorem \ref{th_cohereur} we may assume
that $\cK^\bullet$ arises from a bounded complex $K^\bullet$
of $A$-modules of finite type (details left to the reader).
Then, remark \ref{rem_dual-complex}(i) and corollary
\ref{cor_Ext-loc=glob}(ii) reduce to checking that for every
$A$-module $M$ there exists an open neighborhood $U$ of
$x$ in $X$, such that the $A$-modules
$T^i(M):=H^iR\Hom^\bullet_A(M[0],K^\bullet)$ have support
in $X\setminus U$, for every sufficiently large $i\in\Z$.

We shall argue by induction on the dimension $d$ of the support
of $M$; if $d=0$, there is nothing to show, hence suppose
that $d>0$, and that the assertion is already known for
every $A$-module whose support has dimension strictly
less than $d$. By \cite[6.4]{Mat} we may further reduce
to the case where $M=A/\fp$, where $\fp\subset A$ is
a prime ideal. Then, pick a finite system of generators
$t_1,\dots,t_n\in A$ for $\fp$, and set
$L_\bullet:=\bK_\bullet(t_1,\dots,t_n)$ (notation of remark
\ref{rem_koszul-alg}(ii)). Now, $H_0L_\bullet\simeq A/\fp$,
and $H_iL_\bullet=0$ for every $i>n$; moreover, lemma
\ref{lem_koszul-vanish}(ii) says that $H_iL_\bullet$ is an
$A/\fp$-module of finite type, for every $i=0,\dots,n$,
and vanishes for $i\in\Z$ outside this range; therefore, we may
find $f\in A\setminus\fp$ such that $(H_iL_\bullet)_f$ is a free
$(A/\fp)_f$-module of some finite rank $b(i)\in\N$, for every
$i\in\Z$ (\cite[Th.4.10(ii)]{Mat}), and there exists an integer
$m\leq n$ such that $b(m)\neq 0$ and $b(i)=0$ for every $i>m$.
By choosing an injective resolution of $K^\bullet$, we get
a spectral sequence
$$
E^{pq}_2:=R^p\Hom^\bullet_A((H_qL_\bullet)[0],K^\bullet)
\Rightarrow R^{p+q}\Hom^\bullet_A(L_\bullet,K^\bullet)
$$
with differentials $d^{pq}_r:E^{pq}_r\to E^{p-r+1,q+r}_r$ for
every $p,q\in\Z$ and every $r\geq 2$.
Taking into account corollary \ref{cor_Ext-loc=glob}(ii), we
see that
\set\begin{equation}\label{eq_copies_of-row-zero}
(E^{pq}_2)_f\simeq(E_2^{p,0})^{\oplus b(q)}_f
\qquad
\text{for every $p\in\Z$}.
\end{equation}
By the same token, we also get a natural isomorphism
$$
E^{p0}_2(x):=E^{p0}_2\otimes_A\cO_{\!X,x}\simeq
R^p\Hom^\bullet_{\cO_{\!X(x)}}((A/\fp)^\sim[0](x),\cK^\bullet(x))
\qquad
\text{for every $p\in\Z$}.
$$
Since $\cK^\bullet(x)$ has finite injective dimension,
it follows that there exists an integer $r$ such that
$E^{p,0}_2(x)=0$ whenever $p>r$. Since $E^{pq}_2$ is an
$A$-module of finite type for every $p,q\in\Z$, we deduce
that there exists an open neighborhood $U'$ of $x$ in $X$,
such that, for every $p=r+1,\dots,r+n$, the support of
$E^{p,0}_2$ is contained in $X\setminus U'$ (\cite[Th.4.10(i)]{Mat}).
In this case, the support of $(E^{pq}_2)_f$ is also contained in
$X\setminus U'$, for every $(p,q)\in\Z^{\oplus 2}$ with
$p=r+1,\dots,r+n$.

Now, let $\fq\in U'\cap\Spec\,A_f$ be any prime ideal, and suppose
that $(E_2^{p,0})_\fq\neq 0$ for some integer $p>r+n$. From
\eqref{eq_copies_of-row-zero}, we deduce that $(E^{p,m}_2)_\fq\neq 0$
as well, and it is easily seen that
$(E^{p,m}_\infty)_\fq=(E^{p,m}_2)_\fq$. Especially
\set\begin{equation}\label{eq_edge-stays}
R^{p+m}\Hom_A^\bullet(L_\bullet,K^\bullet)\neq 0.
\end{equation}
However, $L_\bullet$ is also a complex of free $A$-modules,
therefore
$$
R\Hom^\bullet_A(L_\bullet,K^\bullet)\simeq
\Hom^\bullet_A(L_\bullet,K^\bullet)
\qquad
\text{in $\sD(A\Mod)$}.
$$
Thus, if we represent $K^\bullet$ via an object of
$\sC^{[a,b]}(A\Mod)$ (for some $a,b\in\Z$ with
$a\leq b$ : notation of \eqref{sec_brutal-truncate}),
then $R\Hom^\bullet_A(L_\bullet,K^\bullet)$ is represented
by an object of $\sC^{[a,b+n]}(A\Mod)$, and consequently,
$R^i\Hom_A^\bullet(L_\bullet,K^\bullet)=0$ for $i>b+n$.
Taking into account \eqref{eq_edge-stays}, we conclude
that $p\leq p+m\leq b+n$, so that
\set\begin{equation}\label{eq_vanish-outside}
T^i(M)_\fq=0
\qquad
\text{whenever $i>n+\max(b,r)$ and $\fq\in U'\cap\Spec\,A_f$}.
\end{equation}
On the other hand, from the short exact of $A$-modules
$$
0\to M\xrightarrow{\ f\ }M\to M/fM\to 0
$$
we see that $T^i(M)/fT^i(M)\subset T^{i+1}(M/fM)$ for
every $i\in\Z$; since the support of $M/fM$ has dimension
$<d$, the inductive assumption tells us that there exist
$N\in\N$ and an open neighborhood $U''$ of $x$ in $X$ such
that the support of $T^{i+1}(M/fM)$ lies in $X\setminus U''$,
for every $i>N$. From this, Nakayama's lemma implies that
$$
T^i(M)_\fq=0
\qquad
\text{for every $\fq\in\Spec\,A/fA$ and every $i>N$}.
$$
Combining with \eqref{eq_vanish-outside}, we conclude
that $U:=U'\cap U''$ will do.
\end{pfclaim}

Since $X$ is quasi-compact, assumption (b), claim
\ref{cl_fin-inj-dim} and proposition
\ref{prop_Stanley-criterion}(ii) imply that the rule :
$\cF^\bullet\mapsto\cD(\cF^\bullet):=
R\cHom^\bullet_{\cO_{\!X}}(\cF^\bullet,\cK^\bullet)$ takes
objects of $\sD^b(\cO_{\!X}\Mod)_\coh$ to objects of the
same category. Then, corollary \ref{cor_Ext-loc=glob}(ii)
and assumption (b) imply that the natural map
$\cF[0](x)\to\cD\circ\cD(\cF[0])(x)$ is an isomorphism
in $\sD^b(\cO_{\!X(x)}\Mod)$ for every $x\in X$ and every
coherent $\cO_{\!X}$-module $\cF$, whence (a).

(a)$\Rightarrow$(b): Let $x\in X$ be any point, and
$\cM$ any coherent $\cO_{\!X(x)}$-module; we may extend
$\cM$ to a coherent $\cO_{\!U}$-module on some affine open
neighborhood $U$ of $x$ in $X$, and then the latter can
be further extended to a coherent $\cO_{\!X}$-module $\cF$
(lemma \ref{lem_extend-cohs}(ii)). According to
corollary \ref{cor_Ext-loc=glob}(ii), the natural map
$\cD(\cF[0])(x)\to R\cHom_{\cO_{\!X(x)}}^\bullet(\cM[0],\cK^\bullet(x))$
is an isomorphism in $\sD^b(\cO_{\!X(x)}\Mod)$, whence (b).
\end{proof}

\begin{corollary}\label{cor_up-with-duals}
Let $f:Y\to X$ be a morphism of finite type between affine
noetherian schemes, and $\omega^\bullet$ a dualizing complex
on $X$. Then $f^!\omega^\bullet$ is a dualizing complex on $Y$.
\end{corollary}
\begin{proof} We may factor $f$ as the composition of
a closed immersion $g:Y\to Z$ followed by a smooth and
affine morphism $h:Z\to X$. In light of lemmata
\ref{lem_transit-dual}(i) and \ref{lem_compat-for_!}(i),
it then suffices to check that $h^!\omega^\bullet$ is
dualizing on $Z$, so we may assume from start that $f$
is smooth, in which case we are reduced to showing that
$f^*\omega^\bullet$ is dualizing on $Y$. Let $y\in Y$
be any point, set $x:=f(y)$, and denote by $f_y:Y(y)\to X(x)$
the morphism induced by $f$; taking into account
proposition \ref{prop_fishy}, we are further reduced
to checking that $f^*_y(\omega^\bullet_{|X(x)})$ is
dualizing on $Y(y)$. However, the latter assertion follows
easily from proposition \ref{prop_Goren-pullback},
remark \ref{rem_Gorenstein}(i) and \cite[Th.28.7]{Mat}.
\end{proof}

\begin{corollary}\label{cor_fishy}
Let $A$ be a noetherian ring, $x\in X:=\Spec\,A$ any point,
and $K^\bullet$ any object of $\sD^b(\cO_{\!X}\Mod)_\coh$.
Suppose that the following holds :
\begin{enumerate}
\alphaenu
\item
Every irreducible reduced closed subscheme $Z\subset X$
containing $x$ admits a dense open subset which is Gorenstein
(see definition {\em\ref{def_Gorenstein}(i)}).
\item
$K^\bullet(x)$ is dualizing on $X(x)$.
\end{enumerate}
Then there exists an open neighborhood $U$ of $x$
in $X$, such that $K^\bullet_{|U}$ is dualizing on $U$.
\end{corollary}
\begin{proof} We begin with the following :

\begin{claim}\label{cl_generic-Gore}
For every reduced and irreducible closed subscheme $Z$
of $X$ containing $x$, there exists a dense open subset
$U_Z$ of $Z$ such that $K^\bullet(z)$ is dualizing on
$X(z)$, for every $z\in U_Z$.
\end{claim}
\begin{pfclaim} Let $i_Z:Z\to X$ be the closed immersion,
and set $i_{Z(x)}:=i_Z\times_XX(x):Z(x)\to X(x)$; in view of
proposition \ref{prop_replace-prop.12.3.5}(ii), we easily
see that there is a natural isomorphism
$$
i_{Z(x)}^!K^\bullet(x)\isom(i_Z^!K^\bullet)(x)
\qquad
\text{in $\sD(\cO_{Z(x)}\Mod)$}.
$$
Hence, due to (b) and lemma \ref{lem_transit-dual},
$(i^!_ZK^\bullet)(x)$ is dualizing on $Z(x)$. Denote by
$\eta_Z$ the generic point of $Z$; it follows that
$(i^!_ZK^\bullet)(\eta_Z)$ is dualizing on $Z(\eta_Z)$
(proposition \ref{prop_fishy}), and therefore it is
isomorphic to $\cO_{\!Z(\eta_Z)}[c]$ in $\sD(\cO_{\!Z(\eta_Z)}\Mod)$,
due to (a).
On the other hand, claim \ref{cl_fin-inj-dim} implies that
there exists an open neighborhood $U'_Z$ of $x$ in $Z$,
such that $(i^!_ZK^\bullet)_{|U'_Z}$ lies in
$\sD^b(\cO_{\!U'_Z}\Mod)_\coh$. It then follows easily that
there exists a smaller open subset $U_Z\subset U'_Z$ of $Z$,
such that $(i^!_ZK^\bullet)_{|U_Z}\simeq\cO_{\!U_Z}[c]$ in
$\sD(\cO_{\!U'_E}\Mod)$ (details left to the reader).
In view of (a), we may then further shrink $U_Z$, and
assume that $(i^!_ZK^\bullet)_{|U_Z}$ is dualizing on $U_Z$.
Now, let $y\in U_Z$ be any point, denote by
$i_1:\Spec\,\kappa(y)\to Z(y)$ and $i_2:Z(y)\to X(y)$
the closed immersions, and set $i_y:=i_2\circ i_1$.
By proposition \ref{prop_fishy}, the complex
$i_2^!K^\bullet(y)\simeq(i^!_ZK^\bullet)(y)$ is dualizing
on $Z(y)$, so $i_y^!K^\bullet(y)\simeq i_1^!(i_2^!K^\bullet(y))$
is dualizing on $\Spec\,\kappa(y)$ (lemma
\ref{lem_transit-dual}(i)), and consequently
$K^\bullet(y)$ is dualizing on $X(y)$ (proposition
\ref{prop_Stanley-criterion}(i)). The claim follows.
\end{pfclaim}

Denote by $\cZ$ the set of all reduced and irreducible
closed subschemes of $X$ containing $x$, and for every
$Z\in\cZ$, let $U_Z$ be the largest open subset of $Z$
such that $K^\bullet(y)$ is dualizing on $X(y)$ for every
$y\in U_Z$. The subset $U:=\bigcup_{Z\in\cZ}U_Z$ is
ind-constructible in $X$, and it is clearly generizing
(proposition \ref{prop_fishy}). Thus, $U$ is open in $X$
(\cite[Ch.IV, Th.1.10.1]{EGAIV}), and to conclude the
proof, it suffices -- again, by virtue of proposition
\ref{prop_fishy} -- to check that $x\in U$.
Suppose that the latter fails, and let $Z\subset X$ be
the closure of $\{x\}$ in $X$, endowed with its reduced
subscheme structure; by assumption, $Z\cap U=\emptyset$;
on the other hand, $Z\in\cZ$, so that $U_Z\neq\emptyset$,
by claim \ref{cl_generic-Gore}. The contradiction proves
the assertion.
\end{proof}

\begin{remark} Let $A$ be a noetherian ring, and set
$X:=\Spec\,A$.

(i)\ \
If $X$ admits a dualizing complex, then every irreducible
reduced closed subscheme $Z\subset X$ admits a dense open
subset $U$ which is Gorenstein. Indeed, let $\eta_Z$ be the
generic point of $Z$; we know that $Z$ admits a dualizing
complex $\omega_Z^\bullet$ (lemma \ref{lem_transit-dual}(i)),
hence $\omega_Z^\bullet(\eta_Z)$ is dualizing on $Z(\eta_Z)$
(proposition \ref{prop_fishy}). But then there exists an
isomorphism $\cO_{\!Z(\eta_Z)}[c]\isom\omega^\bullet(\eta_Z)$
in $\sD^b(\cO_{\!Z(\eta_Z)}\Mod)$ for some $c\in\Z$
(proposition \ref{prop_unique-dual}), and any such
isomorphism extends to an isomorphism
$\cO_{\!U}[c]\isom(\omega_Z^\bullet)_{|U}$ on some open
subset $U\subset Z$; since $(\omega^\bullet_Z)_{|U}$ is
dualizing on $U$ (lemma \ref{lem_transit-dual}(ii)), the
assertion follows. 

(ii)\ \
The observation in (i) shows that the converse of
corollary \ref{cor_fishy} holds as well : if $x\in X$
is any point, and $U$ any open neighborhood of $x$ in
$X$ such that $U$ admits a dualizing complex, then
conditions (a) and (b) of corollary \ref{cor_fishy}
holds for $x$.

(iii)\ \
If $A$ is local and $X$ admits a dualizing complex,
then the {\em formal fibres} of $A$ are Gorenstein,
{\em i.e.}, if we denote $A^\wedge$ the completion of
$A$ and set $X^\wedge:=\Spec\,A^\wedge$, then
$X^\wedge\times_X\Spec\,\kappa(x)$ is a Gorenstein
scheme, for every $x\in X$.
To check this assertion, notice that $A^\wedge\otimes_AA/\fp$
is the completion of $A/\fp$, and $\Spec\,A/\fp$ admits
a dualizing complex, if $X$ does (lemma \ref{lem_transit-dual}(i));
hence we may replace $A$ by $A/\fp$ and assume from start
that $A$ is a local noetherian domain, and $x$ the
generic point of $X$. Let $\omega^\bullet$ a dualizing
complex on $X$, and $f:X^\wedge\to X$ the natural morphism;
by proposition \ref{prop_Goren-pullback}, the complex
$f^*\omega^\bullet$ is dualizing on $X^\wedge$. On the
other hand, by (i) there exists a non-empty open subset
$U\subset X$ which is Gorenstein, and therefore
$\omega^\bullet_{|U}\simeq\cO_{\!U}[c]$ for some $c\in\Z$.
Due to lemma \ref{lem_transit-dual}(ii), it follows that
$f^{-1}U$ is Gorenstein as well, and then the same holds
for $f^{-1}(x)$, by proposition \ref{prop_fishy}.

(iv)\ \
If $A$ is local and complete, then $X$ admits a dualizing
complex. Indeed, in this case $A$ is a quotient of a
regular local ring (\cite[Th.29.4(ii)]{Mat}), so the
assertion follows from lemma \ref{lem_transit-dual}(i)
and example \ref{ex_reg-is-Gorenstein}.
\end{remark}

\begin{corollary}\label{cor_dual-catenary}
Let $X$ be any noetherian scheme, $\omega^\bullet$ a
dualizing complex on $X$. We have :
\begin{enumerate}
\item
For every $x\in X$ there exists a unique $c\in\Z$ such that
$$
J(x):=R^c\Gamma_{\{x\}}\omega^\bullet_{|X(x)}\neq 0.
$$
\item
For every $x\in X$, the $A$-module $J(x)$ is the injective
hull of the residue field $\kappa(x)$.
\end{enumerate}
\end{corollary}
\begin{proof} In view of proposition \ref{prop_fishy},
we may assume that $X$ is local and $x$ is its closed
point; say that $X=\Spec\,A$ for some local noetherian
ring $A$, in which case $\omega^\bullet$ is the
$\cO_{\!X}$-module arising from a bounded complex
$M^\bullet$ of $A$-modules of finite type. We notice :

\begin{claim}\label{cl_noether-case}
Let $\fm\subset A$ be the maximal ideal.
There exists a unique $c\in\Z$ such that
\begin{enumerate}
\item
$T^i(N):=R^i\Hom_A^\bullet(N,M^\bullet)=0$
for every $A$-module $N$ of finite type supported at
$\{x\}$ and every $i\neq c$.
\item
The map $T^c(A/\fm^n)\to T^c(A/\fm^{n+1})$ induced by the projection
$A/\fm^{n+1}\to A/\fm^n$ is injective, for every $n\in\N$.
\end{enumerate}
\end{claim}
\begin{pfclaim}(i): By proposition \ref{prop_Stanley-criterion}(i),
we know already that there exists $c\in\N$ such that $T^i(A/\fm)=0$
if and only if $i\neq c$. We must therefore check that the
claim holds for this value of $c$. Now, let $N$ be any
$A$-module of finite type supported at $\{x\}$; we may
find a finite filtration
$N_0:=0\subset N_1\subset\cdots\subset N_k$ consisting of
$A$-submodules, such that $N_{j+1}/N_j\simeq A/\fm$ for
every $i=0,\dots,k-1$. We deduce short exact sequences
$$
T^i(A/\fm)\to T^i(N_{j+1})\to T^i(N_j)
\qquad
\text{for every $i\in\N$ and every $j=0,\dots,k-1$}
$$
whence the contention, by a simple induction on $j$.

(ii) is similar : arguing as in the foregoing, we get
an exact sequence
$$
T^{c-1}(\fm^{n+1}/\fm^n)\to T^c(A/\fm^n)\to T^c(A/\fm^{n+1})
$$
whence the assertion.
\end{pfclaim}

On the other hand, corollary \ref{cor_depth-Kosz} yields
a natural identification :
\set\begin{equation}\label{eq_almost-done}
J(x)\isom\colim_{n\in\N} T^c(A/\fm^n)
\end{equation}
and then assertion (i) follows proposition
\ref{prop_Stanley-criterion}(i) and claim \ref{cl_noether-case}.

(ii): Claim \ref{cl_noether-case} also implies that $T^c$
is an exact functor on the full subcategory of $A\Mod$ whose
objects are the $A$-modules supported at $\{x\}$. By
proposition \ref{prop_SGA2-for-coherents}, remark
\ref{rem_SGA2-for-coherents}(i,ii) and \eqref{eq_almost-done}
we deduce already that $J(x)$ is an injective $A$-module.
Next, let $\kappa(x)^\sim$ be the quasi-coherent
$\cO_{\!X}$-module arising from $\kappa(x)$; we have
natural identifications
$$
\begin{aligned}
\Hom_A(\kappa(x),J(x))[c]\isom\, &
R\Hom^\bullet_{\cO_{\!X}}(\kappa(x)^\sim[0],
R\underline\Gamma_{\{x\}}\omega^\bullet)
& \quad \text{(by corollary \ref{cor_Ext-loc=glob}(i))} \\
\isom\, &
R\Gamma_{\!\{x\}}\circ
R\cHom^\bullet_{\cO_{\!X}}(\kappa(x)^\sim[0],\omega^\bullet)
& \quad \text{(by lemma \ref{lem_adj-Gamma_Z}(iii))} \\
\isom\, &
R\Hom^\bullet_{\cO_{\!X}}(\kappa(x)^\sim[0],\omega^\bullet)
& \quad \text{(as $\Supp\,\kappa(x)=\{x\}$)} \\
\isom\, &
R\Hom^\bullet_A(\kappa(x)[0],M^\bullet)
& \quad \text{(by corollary \ref{cor_Ext-loc=glob}(i))} \\
\isom\, & \kappa(x)[c]
& \quad \text{(by proposition \ref{prop_Stanley-criterion}(i))}.
\end{aligned}
$$
Then the assertion follows immediately from theorem
\ref{th_classify-injectives}.
\end{proof}

\sset\subsubsection{}\label{subsec_dual-catenary}
Let $X$ be a noetherian scheme that admits a dualizing
complex $\omega^\bullet$. According to corollary
\ref{cor_dual-catenary}, we get a well defined function
$$
c_X:|X|\to\Z
$$
by assigning to each $x\in X$ the unique $c_X(x)\in\Z$
such that $R^{c_X(x)}\Gamma_{\{x\}}\omega^\bullet_{|X(x)}$
does not vanish. With this notation, we may state :

\begin{lemma}\label{lem_dual-catenary}
In the situation of \eqref{subsec_dual-catenary},
let $x,y\in X$ be any two points, such that $x$ is an
immediate specialization of $y$. Then
$$
c_X(x)=c_X(y)+1.
$$
\end{lemma}
\begin{proof} By proposition \ref{prop_fishy} we may
replace $X$ by $X(x)$, and $\omega$ by $\omega_{|X(x)}$,
and assume from start that $X$ is local, and $x$ is its
closed point. Let $Z\subset X$ be the topological closure
of $\{y\}$, endowed with its reduced closed subscheme
structure. If we denote $i:Z\to X$ the corresponding
closed immersion, then $i^!\omega$ is dualizing on $Z$
(lemma \ref{lem_transit-dual}(i)). Fix $z\in Z$, let
$i_z:Z(z)\to X(z)$ be the induced closed immersion, and
set $c:=c_X(z)$; we get natural isomorphisms :
$$
\begin{aligned}
R\Gamma_{\!\{z\}}(i^!\omega^\bullet)_{|Z(z)}\isom\, &
R\Gamma_{\!\{z\}}\circ
R\cHom^\bullet_{\cO_{X(z)}}(i_{z,*}\cO_{Z(y)},\omega^\bullet_{|X(z)}) &
\qquad\text{(by proposition \ref{prop_replace-prop.12.3.5}(ii))} \\
\isom\, &
R\Hom^\bullet_{\cO_{\!X(z)}}(i_{z,*}\cO_{Z(y)},
R\underline\Gamma_{\!\{z\}}\omega^\bullet_{|X(z)}) &
\qquad\text{(by lemma \ref{lem_adj-Gamma_Z}(iii))} \\
\isom\, &
\Hom^\bullet_{\cO_{\!X(z)}}(i_{z,*}\cO_{Z(y)},
R^{c}\underline\Gamma_{\!\{z\}}\omega^\bullet_{|X(z)})[-c] &
\qquad\text{(by corollary \ref{cor_dual-catenary})}
\end{aligned}
$$
which shows that
$$
c_Z(z)=c_X(z)
$$
provided $c_Z$ is defined via $i^!\omega^\bullet$. Hence,
we may replace $X$ by $Z$, and assume additionally that
$X=\Spec\,A$, where $A$ is a one-dimensional local noetherian
ring, and $y$ is the generic point of $X$. In this case,
$\omega^\bullet$ is the complex of $\cO_{\!X}$-module arising
from a bounded complex $M^\bullet$ of $A$-modules of finite
type, and clearly
$$
R^{c_X(y)}\Gamma_{\!\{y\}}\omega^\bullet_{|X(y)}\isom
H^{c_X(y)}(\omega^\bullet_y)\simeq
H^{c_X(y)}(M^\bullet)\otimes_AK
$$
where $K$ is the field of fractions of $A$. To compute
$c_X(x)$, we apply proposition \ref{prop_depth-Kosz}(iii);
to this aim, let $f$ be any non-zero element of the
maximal ideal of $A$; we get
$$
R^i\Gamma_{\!\{x\}}\omega^\bullet\isom
\colim_{n\in\N} H^i(f^n,M^\bullet)
\qquad
\text{for every $i\in\Z$}
$$
(notation of remark \ref{rem_koszul-alg}(ii)). By inspecting the
definitions, we find a ladder with exact rows
$$
{\diagram
H^iM^\bullet \ar[r]^-{f^n} \ddouble & H^iM^\bullet \ar[r] \ar[d]^f
& H^{i+1}(f^n,M^\bullet) \ar[d] \\
H^iM^\bullet \ar[r]^-{f^{n+1}} & H^iM^\bullet \ar[r] &
H^{i+1}(f^{n+1},M^\bullet) 
\enddiagram}
\qquad
\text{for every $i\in\Z$ and $n\in\N$}
$$
whence -- after taking colimits -- an exact sequence
$$
H^iM^\bullet\xrightarrow{\ \alpha_i\ }(H^iM^\bullet)\otimes_AK
\xrightarrow{\ \beta_i\ } R^{i+1}\Gamma_{\!\{x\}}M^\bullet
\qquad
\text{for every $i\in\Z$}
$$
where $\alpha_i$ is induced by the inclusion map $A\to K$.
However, if $(H^iM^\bullet)\otimes_AK\neq 0$, clearly $\alpha_i$
cannot be surjective, so $\Img\,\beta_i\neq 0$, and the
lemma follows.
\end{proof}

\begin{example}\label{ex_compute-c_X}
Let $K$ be any field, and $A$ a $K$-algebra of finite type.
Set $S:=\Spec\,K$, $X:=\Spec\,A$, and denote by $f:X\to S$
the structure morphism. Since $\cO_{\!S}[0]$ is dualizing on
$S$, the complex $\omega^\bullet:=f^!\cO_{\!S}[0]$ is dualizing
on $X$ (corollary \ref{cor_up-with-duals}). In this case, the
function $c_X$ associated with $\omega_X$ as in
\eqref{subsec_dual-catenary} can be determined explicitly as
follows. First, we may find a closed immersion
$i:X\to Y:=\A^n_S$ for some $n$, and if $p:Y\to S$ is the
structure morphism, the proof of lemma \ref{lem_dual-catenary}
shows that $c_X=c_Y\circ i$, where $c_Y$ is the function of
the same type associated with $g^!\cO_{\!S}[0]$. Hence, we may
assume that $X=\A^n_S$, in which case
$\omega^\bullet\simeq\cO_{\!X}[n]$. Then $X$ is regular
(\cite[Th.28.7]{Mat}), so for every $x\in X$ we may compute :
$$
c_X(x)=\depth_{\{x\}}\cO_{\!X(x)}-n=n-\tr.\deg(\kappa(x)/K)-n=
-\tr.\deg(\kappa(x)/K).
$$
\end{example}

\begin{theorem}\label{th_dual-give-catenary}
Let $A$ be a noetherian ring, and suppose that $X:=\Spec\,A$\ \
admits a dualizing complex. Then $A$ is universally catenary.
\end{theorem}
\begin{proof} In light of corollary \ref{cor_up-with-duals},
it suffices to show that $A$ is catenary. But this follows
easily from lemma \ref{lem_dual-catenary} (details left to
the reader).
\end{proof}

\sset\subsubsection{}\label{subsec_residual}
In the situation of \eqref{subsec_dual-catenary}, lemma
\ref{lem_dual-catenary} shows that $c_X$ is a weak codimension
function, with the terminology of \eqref{subsec_supp-and-cod}.
Then corollary \ref{cor_dual-catenary}(i) means precisely
that $\omega^\bullet$ is a $c_X$-Cohen-Macaulay complex
of $\cO_{\!X}$-modules on $X$ (see definition
\ref{def_Cousin}(i)). By theorem \ref{th_Cousin}, we may
then find a Cousin complex of $\cO_{\!X}$-modules
$\cR_\omega^\bullet$, unique up to isomorphism of complexes,
such that $\omega^\bullet$ is isomorphic to
$\cR^\bullet_\omega$ in $\sD^+(\cO_{\!X}\Mod)$. We call
$\cR_\omega^\bullet$ the {\em residual complex} arising from
$\omega^\bullet$. By corollary \ref{cor_dual-catenary}(ii)
we see that
$$
\cR^p_\omega=\bigoplus_{c_X(x)=p}E(x)^\sim
\qquad
\text{for every $p\in\Z$}
$$
where $E(x)$ is an injective hull of the $A$-module
$\kappa(x)$, for every $x\in X$ (notation of
\eqref{sec_various-O-mod}). 

\sset\subsubsection{}
For every coherent scheme $X$, consider the category
$$
\mathsf{Dual}_X
$$
whose objects are the pairs $(U,\omega^\bullet_U)$, where
$U\subset X$ is an open subset, and $\omega^\bullet_U$ is
a dualizing complex on $X$. The morphisms
$(U,\omega^\bullet_U)\to(U',\omega^\bullet_{U'})$ are the pairs
$(j,\beta)$, where $j:U\to U'$ is an inclusion map of open subsets
of $X$, and $\beta:j^*\omega^\bullet_{U'}\isom\omega^\bullet_U$
is an isomorphism in $\sD^b(\cO_{\!U}\Mod)$. Let $X_\Zar$ denote
the full subcategory of $\Sch/X$ whose objects are the (Zariski)
open subsets of $X$; it follows easily from lemma
\ref{lem_transit-dual}(ii), that the forgetful functor
\set\begin{equation}\label{eq_dualiz-fibrat}
\mathsf{Dual}_X\to X_\Zar
\qquad
(U,\omega^\bullet_U)\mapsto U
\end{equation}
is a fibration (see definition \ref{def_pull-bekko}(ii))
With this notation, we have :

\begin{proposition}\label{prop_glue-dualizers}
Every descent datum for the fibration \eqref{eq_dualiz-fibrat}
is effective.
\end{proposition}
\begin{proof} Let
$((U_i,\omega_i^\bullet); \beta_{ij}~|~i,j=1,\dots,n)$ be a
descent datum for the fibration \eqref{eq_dualiz-fibrat}; this
means that each $(U_i,\omega_i^\bullet)$ is an object of
$\mathsf{Dual}_X$, and if let $U_{ij}:=U_i\cap U_j$ for
every $i,j\leq n$, then
$$
\beta_{ij}:\omega^\bullet_{i|U_{ij}}\isom\omega^\bullet_{j|U_{ij}}
$$
are isomorphisms in $\sD^b(\cO_{\!U_{ij}}\Mod)$, fulfilling
a suitable cocycle condition (see \eqref{subsec_descnt-data}).

We shall show, by induction on $k=1,\dots,n$, that there
exists a dualizing complex $\omega^\bullet_{X_k}$ on
$X_k:=U_1\cup\cdots\cup U_k$, such that the descent datum
$((U_i,\omega_i^\bullet); \beta_{ij}~|~i,j=1,\dots,k)$
is isomorphic to the descent datum relative to the family
$(U_i~|~i=1,\dots,k)$ determined by $(X_k,\omega^\bullet_{X_k})$.

For $k=1$, there is nothing to prove. Next, suppose that
$k>1$, and $\omega^\bullet_{X_{k-1}}$ with the sought properties
has already been constructed. If $k-1=n$, we are done; otherwise,
let $V_k:=X_k\cap U_{k+1}$, set $\fU_k:=(U_{i,k+1}~|~i=1,\dots,k)$,
and $C^\bullet:=R\cHom^\bullet_{\cO_{V_k}}
(\omega^\bullet_{X_k|V_k},\omega^\bullet_{k+1|V_k})$. The system
of morphisms $(\beta_{i,k+1}~|~i=1,\dots,k)$ determines a class
$c_k$ in the \v{C}ech cohomology group $H^0(\fU_k,C^\bullet)$.
However, we have a spectral sequence (see lemma \ref{lem_aspher}) :
$$
E_2^{pq}:=H^p(\fU_k,H^qC^\bullet)\Rightarrow
H^{p+q}R\Hom^\bullet_{\cO_{V_k}}
(\omega^\bullet_{X_k|V_k},\omega^\bullet_{k+1|V_k}).
$$
\begin{claim}\label{cl_peverse-gluing}
$E^{pq}_2=0$ for every $q<0$.
\end{claim}
\begin{pfclaim} More precisely, we check that $H^qC^\bullet=0$
for every $q<0$. Indeed, by lemma \ref{lem_transit-dual}(ii),
that both $\omega^\bullet_{X_k|V_k}$ and $\omega^\bullet_{k+1|V_k}$
are dualizing on $V_k$. Moreover, on $U_{i,k+1}$ they restrict
to isomorphic objects of $\sD^b(\cO_{\!U_{i,k+1}}\Mod)$, for
every $i=1,\dots,k$. It follows easily from proposition
\ref{prop_unique-dual} that there exists an invertible
$\cO_{V_k}$-module $\cL$ and an isomorphism
$$
\omega^\bullet_{X_k|V_k}\isom
\omega^\bullet_{k+1|V_k}\otimes_{\cO_{V_k}}\cL
\qquad
\text{in $\sD^b(\cO_{V_k}\Mod)$}.
$$
Therefore, the biduality map induces an isomorphism
$$
R\cHom^\bullet_{\cO_{V_k}}
(\omega^\bullet_{X_k|V_k},\omega^\bullet_{k+1|V_k})\isom\cL[0]
$$
whence the claim.
\end{pfclaim}

Claim \ref{cl_peverse-gluing} implies that $c_k$ corresponds
to an element of $H^0R\Hom_{\cO_{V_k}}^\bullet
(\omega^\bullet_{X_k|V_k},\omega^\bullet_{k+1|V_k})$, and by
construction, it is clear that this global section is an
isomorphism $c_k:\omega^\bullet_{X_k|V_k}\isom\omega^\bullet_{k+1|V_k}$
in $\sD^b(\cO_{V_k}\Mod)$. By definition, $c_k$ is represented
by a diagram of quasi-isomorphisms :
\set\begin{equation}\label{eq_this-is-c_k}
\omega^\bullet_{X_k|V_k}\leftarrow T^\bullet\to\omega^\bullet_{k+1|V_k}
\qquad
\text{in $\sC(\cO_{V_k}\Mod)$}
\end{equation}
for some complex $T^\bullet$. Let $T_!$ be the $\cO_{\!X_{k+1}}$-module
obtained as extension by zero of $T$, and define likewise 
$(\omega^\bullet_{X_k})_!$ and $(\omega^\bullet_{k+1})_!$.
We let $\omega^\bullet_{X_{k+1}}$ be the cone of the map of
complexes
$$
T_!\to(\omega^\bullet_{X_k})_!\oplus(\omega^\bullet_{k+1})_!
$$
deduced from \eqref{eq_this-is-c_k}. An easy inspection shows
that this complex is dualizing on $X_{k+1}$, and it fulfills
the stated condition.
\end{proof}

\begin{remark} The proof of proposition \ref{prop_glue-dualizers}
is a special case of a general technique for ``glueing perverse
sheaves'' developed in \cite{BBD}.
\end{remark}

In the study of duality for derived categories of modules,
the role of reflexive modules is taken by the more general
class of Cohen-Macaulay modules. There is a version of this
theory for noetherian regular schemes, and a relative
variant, for smooth morphisms. Let us begin by recalling
the following :

\begin{definition}\label{def_CM}
Let $(A,\fm_A)$ be a local ring, $M$ an $A$-module.
\begin{enumerate}
\item
Let $\cM$ be the set of all finitely generated $A$-submodules
$M$. The {\em dimension\/} of $M$ is
$$
\dim_AM:=\sup(\dim\,\Supp\,M'~|~M'\in\cM).
$$
\item
Suppose that the open subset $\Spec\,A\setminus\{\fm_A\}$
is quasi-compact, and $\depth_A\,M<+\infty$. If
$\dim_AM=\depth_A\,M$, we say that $M$ is a
{\em Cohen-Macaulay\/} $A$-module. The category
$$
A\CM
$$
is the full subcategory of $A\Mod$ consisting of
all finitely presented Cohen-Macaulay $A$-modules.
For every $n\in\N$, we let $A\CM_n$ be the full
subcategory of $A\CM$ whose objects are the Cohen-Macaulay
$A$-modules of dimension $n$.
\item
If $A$ is a Cohen-Macaulay $A$-module, we say that
$A$ is a {\em Cohen-Macaulay} local ring.
\item
Let $X$ be a scheme, $\cF$ a quasi-coherent $\cO_{\!X}$-module.
We say that $\cF$ is a {\em Cohen-Macaulay} $\cO_{\!X}$-module,
if $\cF_x$ is a Cohen-Macaulay $\cO_{\!X,x}$-module, for every
$x\in X$. We say that $X$ is a {\em Cohen-Macaulay\/} scheme,
if $\cO_{\!X}$ is a Cohen-Macaulay $\cO_{\!X}$-module.
\item
Let $B$ be another local ring, $\phi:A\to B$ a local ring
homomorphism. The category
$$
\phi\CM
$$
of {\em $\phi$-Cohen-Macaulay modules\/} is the full
subcategory of $B\Mod$ whose objects are all the finitely
presented $B$-modules $N$ that are $\phi$-flat, and such
that $N/\fm_AN$ is a Cohen-Macaulay $B$-module. For every
$n\in\N$, we let $\phi\CM_n$ be the full subcategory of
$\phi\CM$ whose objects are the $\phi$-Cohen-Macaulay
modules $N$ with $\dim_B N/\fm_AN=n$.
\item
Let $f:X\to Y$ be a locally finitely presented morphism
of schemes, $\cF$ a quasi-coherent $\cO_{\!X}$-module,
$x\in X$ any point, and set $y:=f(x)$. We say
that $\cF$ is {\em $f$-Cohen-Macaulay\/} at the point $x$,
if $\cF_x$ is a $f^\natural_x$-Cohen-Macaulay module.
We say that $f$ is {\em Cohen-Macaulay\/} at the point
$x$, if $\cO_{\!X}$ is $f$-Cohen-Macaulay at the point
$x$ (cp. \cite[Ch.IV, D\'ef.6.8.1]{EGAIV-2}).
\item
Let $f$ and $\cF$ be as in (vi). We say that $\cF$ is
$f$-Cohen-Macaulay (resp. that $f$ is Cohen-Macaulay)
if $\cF$ (resp. $\cO_{\!X}$) is $f$-Cohen-Macaulay at
every point $x\in X$.
\item
Let $f$ and $\cF$ be as in (vi). The {\em $f$-Cohen-Macaulay
locus\/} of $\cF$ is the subset
$$
CM(f,\cF)\subset X
$$
consisting of all $x\in X$ such that $f$ is
Cohen-Macaulay at $x$. The subset $CM(f,\cO_{\!X})$
is also called the {\em Cohen-Macaulay locus of $f$}
and is denoted briefly $CM(f)$.
\end{enumerate}
\end{definition}

\begin{remark}
Let $(A,\fm_A)$ be any local ring; set $X:=\Spec\,A$, and
for every $A$-module $M$, let $M^\sim$ be the quasi-coherent
$\cO_X$-module arising from $M$. Then, with the notation of
theorem \ref{th_tohoku-vanish}(i), we have
$$
\dim_AM=d_{M^\sim}
\qquad
\text{for every $A$-module $M$}.
$$
For the proof, let $\cM$ be the set of all finitely generated
$A$-submodules of $M$. Clearly $d_{M^\sim}=\sup(d_{M'^\sim}~|~M'\in\cM)$,
so we are reduced to the case where $M$ is finitely generated.
Then $\Supp\,M$ is a closed subset of $X$, and
$\dim_AM=\dim(\Supp\,M)$; moreover, clearly $\Supp\,M=\Supp\,M^\sim$,
and the assertion folows easily.
\end{remark}

\begin{lemma}\label{lem_CM}
Let $f:X\to Y$ be a locally finitely presented morphism
of schemes, and $\cF$ a finitely presented $\cO_{\!X}$-module
such that $\Supp\,\cF=X$.
We have :
\begin{enumerate}
\item
$CM(f,\cF)$ is an open subset of $X$,
\item
The restriction $CM(f)\to Y$ of $f$ is locally equidimensional.
\item
Let $g:X'\to X$ be another morphism of schemes, $x'\in X'$
a point such that $g$ is \'etale at $x'$. Then
$g(x')\in CM(f,\cF)$ if and only if $x'\in CM(f\circ g,g^*\cF)$.
\item
If $x\in CM(f)$, then
$\delta'(x,\cO_{\!X})=\delta'(f(x),\cO_Y)+\dim\cO_{f^{-1}(fx),x}$
(notation of \eqref{subsec_local-depth}).
\end{enumerate}
\end{lemma}
\begin{proof} (i) and (ii) follow easily from
\cite[Ch.IV, Prop.15.4.3]{EGAIV-3} and
\cite[Ch.IV, Prop.2.3.4]{EGAIV-2}.

(iii) follows from corollary \ref{cor_f-flat-invariance}.
Lastly, (iv) is an immediate consequence of corollary
\ref{cor_depth-flat-basechange}.
\end{proof}

\begin{remark}
(i)\ \
The terminology ``locally equidimensional''
follows \cite[Err${}_\mathrm{IV}35$]{EGA4}, which modifies
the terminology ``equidimensional'' of
\cite[Ch.IV, D\'ef.13.3.2]{EGAIV-3}. Hence, a morphism of
schemes $f:X\to Y$ locally of finite type is locally
equidimensional, if it verifies the equivalent conditions
of \cite[Ch.IV, Prop.13.3.1]{EGAIV-3} at every point $x\in X$.
However, the {\em caveat} here is that the said conditions
are {\em not} quite equivalent as stated : they become equivalent
if one omits the words ``contenant $y$'' in line 2 of {\em loc.cit.}
(with unchanged proof). So, the correct definition of ``locally
equidimensional" is in terms of the cited proposition thus amended.

(ii)\ \
A counterexample to \cite[Ch.IV, Prop.13.3.1]{EGAIV-3}
is given by the identity morphism $X\to X$ of a scheme
$X$ with a point $x$ such that the union of the irreducible
components of $X$ passing through $x$ does not contain an
open neighborhood of $x$ in $X$ ({\em e.g.} take $X$ whose
underlying topological space is zero-dimensional and not
discrete); then $\one_X$ satisfies condition b) of
\cite[Ch.IV, Prop.13.3.1]{EGAIV-3}, but does not satisfy
conditions a),$\mathrm{a}'$) and $\mathrm{a}''$).
\end{remark}

\begin{proposition}\label{prop_not-so-triv}
Let $A$ be a regular local ring of dimension $d$, and $M$ a
finitely generated Cohen-Macaulay $A$-module. Then :
\begin{enumerate}
\item
$\Ext^i_A(M,A)=0$ for every $i\neq c:=d-\dim M$.
\item
The $A$-module $\cD(M):=\Ext^c_A(M,A)$ is Cohen-Macaulay.
\item
The natural map $M\to\Ext^c(\cD(M),A)$ is an isomorphism,
and $\Supp\,\cD(M)=\Supp\,M$.
\item
For every $n\in\N$, we have an equivalence of categories :
$$
\cD:A\CM_n\to A\CM^o_n
\qquad
N\mapsto\Ext_A^{d-n}(N,A).
$$
\end{enumerate}
\end{proposition}
\begin{proof}(i): According to \cite[Ch.0, Prop.17.3.4]{EGAIV},
the projective dimension of $M$ equals $c$, so we may find
a minimal free resolution
$0\to L_c\to\cdots\to L_1\to L_0\to M$ for $M$ of length
$c$. Hence, we need only prove the sought vanishing for
every $i<c$. Set $X:=\Spec\,A$, $Z:=\Supp\,M\subset X$.
By virtue of proposition \ref{prop_compute-depth}, it
suffices to show that $\depth_Z\cO_{\!X}\geq c$.
In light of \eqref{eq_for-I-depth}, this comes down to
showing that $\cO_{\!X,z}$ is a local ring of depth
$\geq c$, for every $z\in Z$. The latter holds, since
$\cO_{\!X,z}$ is a regular local ring of dimension
$\geq c$ (\cite[Ch.0, Cor.16.5.12]{EGAIV}).

(ii): From (i) we deduce that
$L_\bullet^\vee:=(L^\vee_0\to\cdots\to L^\vee_n)$, together
with its natural augmentation $L^\vee_n\to\cD(M)$, is a
free resolution of $\cD(M)$; especially, the projective
dimension of $\cD(M)$ is $\leq c$, and therefore
$$
\depth_A\cD(M)\geq\dim M
$$
again by \cite[Ch.0, Prop.17.3.4]{EGAIV}. On the other
hand, it follows easily from \cite[Prop.3.3.10]{We} that
$\Supp\,\cD(M)\subset\Supp\,M$, so $\cD(M)$ is Cohen-Macaulay.

(iii): From the proof of (ii) we see that $R\Hom_A(\cD(M),A)$
is computed by the complex $L^{\vee\vee}_\bullet=L_\bullet$,
whence the first assertion. Invoking again \cite[Prop.3.3.10]{We},
we deduce that $\Supp\,M\subset\Supp\,\cD(M)$; since the
converse inclusion is already known, these two supports
coincide.

(iv) follows straightforwardly from (ii) and (iii).
\end{proof}

\begin{corollary}\label{cor_not-so-triv}
Let $X$ be a regular noetherian scheme, and $j:Y\to X$ a
closed immersion. We have :
\begin{enumerate}
\item
$Y$ admits a dualizing complex $\omega^\bullet_Y$.
\item
If $Y$ is a Cohen-Macaulay scheme, then we may find
a dualizing complex $\omega^\bullet_Y$ that is concentrated
in degree $0$. Moreover, $H^0(\omega^\bullet_Y)$ is a
Cohen-Macaulay $\cO_Y$-module.
\end{enumerate}
\end{corollary}
\begin{proof} Indeed, lemma \ref{lem_transit-dual}(i) shows
that the complex
$\omega^\bullet_Y:=j^*R\cHom_{\cO_{\!X}}(j_*\cO_Y,\cO_{\!X})$
will do. According to proposition \ref{prop_not-so-triv}(i,ii),
this complex fulfills the condition of (ii), up to a suitable
shift.
\end{proof}

\sset\subsubsection{}\label{subsec_get-that}
Let $f:X\to S$ be a smooth quasi-compact morphism of schemes,
$\cF$ a finitely presented quasi-coherent $\cO_{\!X}$-module,
$x\in X$ any point, and $s:=f(x)$. Set
$$
A:=\cO_{\!S,s} \qquad B:=\cO_{\!X,x} \qquad F:=\cF_x
\qquad
B_0:=B\otimes_A\kappa(s) \qquad F_0:=F\otimes_A\kappa(s).
$$

\begin{theorem}\label{th_CM-duality}
In the situation of \eqref{subsec_get-that}, suppose that
$F$ is a $f^\natural_x$-Cohen-Macaulay module. Then we may
find an open neighborhood $U\subset X$ of $x$ in $X$ such
that the following holds:
\begin{enumerate}
\item
There exists a finite resolution of the $\cO_{\!U}$-module $\cF_{|U}$,
of length $n:=\dim B_0-\dim_BF_0$
$$
\Sigma_\bullet
\quad : \quad
0\to\cL_n\to\cL_{n-1}\to\cdots\to\cL_0\to\cF_{|U}\to 0
$$
consisting of free $\cO_{\!U}$-modules of finite rank. Moreover,
$\Sigma_\bullet$ is {\em universally $\cO_{\!S}$-exact}, {\em i.e.}
for every coherent $\cO_{\!S}$-module $\cG$, the complex
$\Sigma_\bullet\otimes_{\cO_{\!X}}f^*\cG$ is still exact.
\item
The complex
$K^\bullet:=R\cHom^\bullet_{\cO_{\!U}}(\cF_{|U},\cO_{\!U})$
is concentrated in degree $n$.
\item
The natural map
$\cF_{|U}\to R\cHom_{\cO_{\!U}}(K^\bullet,\cO_{\!U})$ is an
isomorphism in $\sD(\cO_{\!U}\Mod)$.
\item
The $B$-module $G:=H^nK^\bullet_x$ is $f^\natural_x$-Cohen-Macaulay,
and\/ $\Supp\,G=\Supp\,F$.
\end{enumerate}
\end{theorem}
\begin{proof}(i): According to proposition \ref{prop_resol-exist},
the $B$-module $F$ admits a minimal free resolution
$$
\Sigma_{x,\bullet}
\quad : \quad
\cdots\to L_2\xrightarrow{\ d_2\ } L_1\xrightarrow{\ d_1\ } L_0
\xrightarrow{\ \eps\ }F\to 0
$$
that is universally $A$-exact; especially, $d_i(L_i)$ is a
flat $A$-module, for every $i\in\N$. It also follows that
$\Sigma_{x,\bullet}\otimes_A\kappa(s)$ is a minimal free
resolution of the $B_0$-module $F_0$. Since the latter
is Cohen-Macaulay, and $B_0$ is a regular local ring
(\cite[Ch.IV, Th.17.5.1]{EGA4}), the projective dimension
of the $B_0$-module $F_0$ equals $n$
(\cite[Ch.0, Prop.17.3.4]{EGAIV}), therefore
$d_n(L_n)\otimes_A\kappa(s)$ is a free $B_0$-module, so
$d_n(L_n)$ is a flat $B$-module (lemma \ref{lem_Tor-crit});
then we deduce that it is actually a free $B$-module, as it
is finitely presented. Since $\Sigma_{x,\bullet}$ is minimal,
we conclude that $L_{i}=0$ for every $i>n$.
We may now extend $\Sigma_{x,\bullet}$ to a finite resolution
$\Sigma_\bullet$ of $\cF_{|U}$ by free $\cO_{\!U}$-modules on
some open neighborhood $U$ of $x$, and after replacing $U$
by a smaller neighborhood of $U$, we may assume that $\cF_{|U}$
is $f_{|U}$-flat (\cite[Ch.IV, Th.11.3.1]{EGAIV-3}), hence
$\Sigma_\bullet$ is universally $\cO_{\!S}$-exact, as stated.

(ii): Let $L_\bullet:=(L_n\to\cdots\to L_0)$ be the complex
obtained after omitting $F$ from the resolution
$\Sigma_{x,\bullet}$ (where $L_0$ is placed in degree $0$).
Then $K^\bullet$ is isomorphic to
$L^\vee_\bullet:=\Hom_B(L_\bullet,B)$. On the other hand,
the universal exactness property of $\Sigma_{x,\bullet}$
implies that
$$
\Hom_{B_0}(L_\bullet\otimes_A\kappa(s),B_0)=
L^\vee_\bullet\otimes_A\kappa(s)
$$
computes $R\Hom^\bullet_{B_0}(F_0,B_0)$. In view of proposition
\ref{prop_not-so-triv}(i), we have $\Ext^i_{B_0}(F_0,B_0)=0$
for every $i\neq n$. From this, a repeated application of
\cite[Ch.IV, Prop.11.3.7]{EGAIV-3} shows that
$L^\vee_\bullet$ is concentrated in degree $n$ as well, and
after shrinking $U$, we may assume that (ii) holds (details
left to the reader).

(iii): Let $\cL_\bullet:=(\cL_n\to\cdots\to\cL_0)$ be the
complex obtained by omitting $\cF_{|U}$ form the resolution
$\Sigma_\bullet$; the proof of (iii) shows that
$R\cHom_{\cO_{\!U}}(K^\bullet,\cO_{\!U})$ is computed
by $\cL^{\vee\vee}_\bullet=\cL_\bullet$, whence the
contention.

(iv): By the same token, \cite[Ch.IV, Prop.11.3.7]{EGAIV-3}
implies that $\Coker\,d^\vee_i$ is a flat $A$-module,
for every $i=1,\dots,n$. Especially, $G$ is a flat $A$-module,
and the complex $L^\vee_\bullet[-n]$, with its natural
augmentation $L^\vee_n\to G$, is a universally $A$-exact
and free resolution of the $B$-module $G$. Especially,
the projective dimension of $G$ is $\leq n$, therefore
\set\begin{equation}\label{eq_uof}
\depth_BG\geq\dim_BF_0
\end{equation}
by \cite[Ch.0, Prop.17.3.4]{EGAIV}. From (iii) we also see
that the induced map
$$
F_\fp\to\Ext^n_{B_\fp}(G_\fp,B_\fp)
$$
is an isomorphism, for every prime ideal $\fp\subset B$;
therefore $\Supp\,F\subset\Supp\,G$. Symmetrically, the
same argument yields $\Supp\,G\subset\Supp\,F$. Hence the 
supports of $F$ and $G$ agree. Lastly, combining with
\eqref{eq_uof} we see that $G$ is Cohen-Macaulay.
\end{proof}

\sset\subsubsection{}
Keep the notation of \eqref{subsec_get-that}, and let
$\phi:=f^\natural_x:A\to B$. Theorem
\ref{th_CM-duality}(iii,iv) implies that, for every
$n\in\N$, the functor
$$
\cD_\phi:\phi\CM_n\to\phi\CM^o_n
\qquad
M\mapsto\Ext^{\dim B_0-n}_B(M,B)
$$
is an equivalence, and the natural map
$M\to\cD_\phi\circ\cD_\phi(M)$ is an isomorphism, for every
$\phi$-Cohen-Macaulay module $M$. We shall see later also
a relative variant of corollary \ref{cor_not-so-triv}, in a
more special situation (see proposition \ref{prop_duality-smooth}).

\subsection{Schemes over a valuation ring}\label{sec_sch-val-rings}
Throughout this section, $(K,|\cdot|)$ is a valued field, whose
valuation ring (resp. maximal ideal, resp. residue field, resp. value
group) shall be denoted $K^+$ (resp. $\fm_K$, resp. $\kappa$, resp.
$\Gamma$). Also, we let
$$
S:=\Spec\,K^+
\qquad\text{and}\qquad
S_{\!/b}:=\Spec\,K^+/bK^+
\qquad
\text{for every $b\in\fm_K$}
$$
(so $S_{\!/0}=S$) and we denote by $s:=\Spec\,\kappa$ (resp.
by $\eta:=\Spec\,K$) the closed (resp. generic) point of $S$.
More generally, for every $S$-scheme $X$ we let as well
$$
X_{\!/b}:=X\times_SS_{\!/b}
\qquad
\text{for every $b\in\fm_K$}.
$$
A basic fact, that follows immediately from corollary
\ref{cor_coherence}, is that every finitely presented
$S$-scheme is coherent (this can also be deduced easily
from \cite[Part I, Th.3.4.6]{Gr-Ra}).

\begin{proposition}\label{prop_hdim-finite}
Let $f:X\to S$ be a flat and finitely presented morphism,
$\cF$ a coherent $\cO_{\!X}$-module, $x\in X$ any point,
$y:=f(x)$, and suppose that $f^{-1}(y)$ is a regular
scheme (this holds {\em e.g.} if $f$ is a smooth morphism).
Let also $n:=\dim\cO_{f^{-1}(y),x}$. Then:
\begin{enumerate}
\item
$\hdim_{\cO_{\!X,x}}\cF_x\leq n+1$ (here we set
$\hdim_{\cO_{\!X,x}}0:=0$).
\item
If $\cF$ is $f$-flat at the point $x$, then
$\hdim_{\cO_{\!X,x}}\cF_x\leq n$.
\item
If $\cF$ is reflexive at the point $x$, then
$\hdim_{\cO_{\!X,x}}\cF_x\leq\max(0,n-1)$.
\item
If $f$ has regular fibres, and $\cF$ is reflexive and generically
invertible (see \eqref{subsec_rk-trivias}), then $\cF$ is invertible.
\item
All the fibres of the induced morphism $f_x:X(x)\to S(y)$ are
regular and irreducible.
\end{enumerate}
\end{proposition}
\begin{proof} To start out, we may assume that $X$ is affine,
and then it follows easily from lemma \ref{lem_rflx-on-lim}(ii.a)
that $\cO_{\!X}$ is coherent.

(ii): Since $\cO_{\!X}$ is coherent, we can find a
possibly infinite resolution
$$
\cdots\to E_2\to E_1\to E_0\to\cF_x\to 0
$$
by free $\cO_{\!X,x}$-modules $E_i$ ($i\in\N$) of finite rank;
set $E_{-1}:=\cF_x$ and $L:=\Img(E_n\to E_{n-1})$. It suffices
to show that the $\cO_{\!X,x}$-module $L$ is free. For every
$\cO_{\!X,x}$-module $M$ we shall let
$M(y):=M\otimes_{\cO_{\!X,x}}\cO_{\!f^{-1}(y),x}$.
Since $L$ and $\cF_x$ are torsion-free, hence flat $K^+$-modules,
the induced sequence of $\cO_{\!f^{-1}(y),x}$-modules:
$$
0\to L(y)\to E_{n-1}(y)\to \cdots E_1(y)\to E_0(y)\to\cF_x(y)\to 0
$$
is exact; since $f^{-1}(y)$ is a regular scheme, $L(y)$ is a
free $\cO_{\!f^{-1}(y),x}$-module of finite rank. Since $L$ is
also flat as a $K^+$-module, \cite[Ch.IV, Prop.11.3.7]{EGAIV-3}
and Nakayama's lemma show that any set of elements of $L$
lifting a basis of $L(y)$, is a free basis of $L$.

(i): Locally on $X$ we can find an epimorphism
$\cO^{\oplus n}_{\!X}\to\cF$, whose kernel $\cG$ is again a coherent
$\cO_{\!X}$-module, since $\cO_{\!X}$ is coherent. Clearly it suffices
to show that $\cG$ admits, locally on $X$, a finite free resolution
of length $\leq n$, which holds by (ii), since $\cG$ is $f$-flat.

(iii): Suppose $\cF$ is reflexive at $x$; by remark \ref{rem_co-repres}
we can find a left exact sequence
$$
0\to\cF_x\to\cO^{\oplus m}_{\!X,x}\xrightarrow{\alpha}
\cO^{\oplus n}_{\!X,x}.
$$
It follows that $\hdim_{\cO_{\!X,x}}\cF_x=
\max(0,\hdim_{\cO_{\!X,x}}(\Coker\,\alpha)-2)\leq\max(0,n-1)$, by (i).

(iv): Suppose that $\cF$ is generically invertible, let
$j:U\to X$ be the maximal open immersion such that $j^*\cF$
is an invertible $\cO_{\!U}$-module, and set $Z:=X\!\setminus\! U$.
Under the current assumptions, $X$ is reduced; it follows easily
that $U$ is the set where $\rk\,\cF =1$, and that $j$ is quasi-compact
(since the rank function is constructible : see \eqref{subsec_rk-trivias}).
It follows from (iii) that $Z\cap f^{-1}(y)$ has codimension
$\geq 2$ in $f^{-1}(y)$, for every $y\in S$. Hence the conditions
of corollary \ref{cor_depth-cons} are fulfilled, and furthermore $\cF[0]$
is a perfect complex by (ii), so the invertible $\cO_{\!X}$-module
$\sdet\,\cF$ is well defined (lemma \ref{lem_det-fctr}).
Then (iv) follows from the natural isomorphisms :
$$
\cF\isom j_*j^*\cF\isom j_*\sdet(j^*\cF[0])\isom
j_*j^*\sdet\,\cF[0]\stackrel{\sim}{\leftarrow}\sdet\,\cF[0].
$$

(v): Let $y'\in S(y)$ be any point; we need to check that
$f_x^{-1}(y')$ is an irreducible regular scheme, and after
replacing $K^+$ by a quotient, we may assume that $y'$ is
the generic point of $S$; after further replacing $K^+$ by
a localization, we may also assume that $y$ is the closed
point of $S$. Set $A:=\cO_{X,x}$; by assumption, $f_x^{-1}(y)$
is a regular scheme, hence $\fp:=\fm_KA$ is a prime ideal
of $A$, so that $A_\fp$ is a valuation ring (proposition
\ref{prop_integral-fp-ext}(ii)), and especially, it is a
domain. Moreover, \cite[Ch.IV, Prop.11.3.7]{EGAIV-3} implies
that for every $h\in A\setminus\fp$, scalar multiplication
by $h$ is injective on $A$; thus, the localization map
$A\to A_\fp$ is injective as well, so $A$ is a domain. This
already proves that $f_x^{-1}(y')=\Spec\,K\otimes_{K^+}A$
is irreducible. Next, let $\fq\subset A$ be any prime ideal
with $\fq\cap K^+=\{0\}$, and $M$ any $A_\fq$-module of
finite type; pick a finitely presented $A$-module $N$ with
an isomorphism $N_\fq\isom M$. By (i), there exists a finite
resolution $L_\bullet$ of $N$, consisting of free $A$-modules
of finite rank; then $A_\fq\otimes_AL_\bullet$ is a finite
resolution of $M$ by free $A_\fq$-modules of finite rank.
By Serre's theorem \cite[Th.19.2]{Mat}, this shows that
$A_\fq$ is regular, as required.
\end{proof}

\begin{proposition}\label{prop_Mu}
Let $X$ be a quasi-compact and quasi-separated scheme,
$n>0$ an integer, and denote by $\pi:\P^n_{\!X}\to X$
the natural projection. Then the map
$$
H^0(X,\Z)\oplus\Pic\,X\to\Pic\,\P^n_{\!\!X} \qquad
(r,\cL)\mapsto\cO_{\P^n_{\!\!X}}(r)\otimes\pi^*\cL
$$
is an isomorphism of abelian groups.
\end{proposition}
\begin{proof} To begin with, let us recall the following well known :

\begin{claim}\label{cl_EGAIII-2}
Let $f:Y\to T$ be a proper morphism of noetherian schemes,
$\cF$ a coherent $f$-flat $\cO_Y$-module, $t\in T$ a point,
$p\in\N$ an integer, and denote $i_t:f^{-1}(t)\to Y$ the
natural immersion. Suppose that $H^{p+1}(f^{-1}(t),i_t^*\cF)=0$.
Then the natural map
$$
(R^pf_*\cF)_t\to H^p(f^{-1}(t),i_t^*\cF)
$$
is surjective.
\end{claim}
\begin{pfclaim} In light of \cite[Ch.III, Prop.1.4.15]{EGAIII},
we easily reduce to the case where $T$ is a local scheme, say
$T=\Spec\,A$ for some noetherian local ring $A$, and $t$ is the
closed point of $T$. Set $k:=\kappa(t)$, and for every $q\in\N$,
consider the functor
$$
F_{-q}:A\Mod\to A\Mod
\qquad
M\mapsto H^0(T,R^pf_*(f^*M^\sim\otimes_{\cO_Y}\cF))
$$
(notation of \eqref{sec_various-O-mod}). Since $\cF$
is $f$-flat, the system $(F_{-q}~|~q\in\N)$ defines a
homological functor. The assertion is that $F_{-p-1}(k)=0$,
in which case \cite[Ch.III, Cor.7.5.3]{EGAIII-2} (together
with \cite[Ch.III, Th.3.2.1 and Th.4.1.5]{EGAIII}) says that
$F_{-p-1}(M)=0$ for every $A$-module $M$. Hence, $F_{-p-1}$
is trivially an exact functor, and it follows that the natural
map $F_{-p}(A)\to F_{-p}(k)$ is surjective
(\cite[Ch.III, Prop.7.5.4]{EGAIII-2}), which is the claim.
\end{pfclaim}

Now, the injectivity of the stated map is clear (details left
to the reader). For the surjectivity, let $\cG$ be an invertible
$\cO_{\P^n_{\!X}}$-module, write $X$ as the limit of a filtered
system $(X_\lambda~|~\lambda\in\Lambda)$ of noetherian schemes,
and for every $\lambda\in\Lambda$, let
$p_\lambda:\P^n_{\!X}\to\P^n_{\!X_\lambda}$ be the induced
morphism; for some $\lambda\in\Lambda$, we may find an
invertible $\cO_{\P^n_{\!X_{\!\lambda}}}$-module $\cF_\lambda$
with an isomorphism $\cG\isom p^*_\lambda\cG_\lambda$ of
$\cO_{\P^n_{\!X}}$-modules. Clearly, it then suffices to check
the sought surjectivity for the scheme $X_\lambda$; we may
therefore replace $X$ by $X_\lambda$, and assume from start
that $X$ is noetherian.

Next, let $x\in X$ be any point, and
$i_x:\P^n_{\!\kappa(x)}\to\P^n_{\!X}$ the natural immersion;
we have $i^*_x\cG\simeq\cO_{\P^n_{\!\kappa(x)}}(r_x)$ for some
$r_x\in\Z$. Set $\cF:=\cG(-r_x)$, and notice that
$H^1(\P^n_{\!\kappa(x)},i^*_x\cF)=0$; moreover,
$H_x:=H^0(\P^n_{\!\kappa(x)},i^*_x\cF)$ is a one-dimensional
$\kappa(x)$-vector space. From claim \ref{cl_EGAIII-2},
it follows that there exist an open neighborhood $U$ of
$x$ in $X$, and a section $s\in H^0(\pi^{-1}U,\cF)$ mapping
to a generator of $H_x$. The section $s$ defines a map of
$\cO_{\P^n_{\!U}}$-modules
$\phi:\cO_{\P^n_{\!U}}(r_x)\to\cG_{|\pi^{-1}U}$. Furthermore,
since $i^*_x\cF\simeq\cO_{\P^n_{\kappa(x)}}$, it is easily
seen that $\phi$ restricts to an isomorphism on some open
subset of the form $\pi^{-1}U'$, where $U'\subset U$ is an
open neighborhood of $x$ in $X$.

Summing up, since $x\in X$ is arbitrary, and since $X$ is
quasi-compact, we may find a finite covering
$X=U_1\cup\cdots\cup U_d$ consisting of open subsets, and
for every $i=1,\dots,d$ an integer $r_i$ and an isomorphism
$\phi_i:\cO_{\P^n_{\!U_i}}(r_i)\isom\cG_{|\pi^{-1}U_i}$ of
$\cO_{\P^n_{\!U_i}}$-modules. Then, for every $i,j\leq d$
such that $U_{ij}:=U_i\cap U_j\neq\emptyset$, the composition
of the restriction of $\phi_i$ followed by the restriction of
$\phi_j^{-1}$, is an isomorphism
$$
\phi_{ij}:\cO_{\P^n_{U_{ij}}}(r_i)\isom\cO_{\P^n_{U_{ij}}}(r_j)
$$
of $\cO_{\P^n_{\!U_{ij}}}$-modules. This already implies that
$r_i=r_j$ whenever $U_{ij}\neq\emptyset$, hence the rule :
$x\mapsto r_i$ if $x\in U_i$ yields a well defined continuous
map $r:X\to\Z$. Lastly, $\phi_{ij}$ is the scalar multiplication
by a section
$$
u_{ij}\in H^0(\P^n_{\!U_{ij}},\cO^\times_{\P^n_{U_{ij}}})=
H^0(U_{ij},\cO^\times_{\!U_{ij}})
\qquad
\text{for every $i,j=1,\dots,d$}.
$$
The system $(u_{ij}~|~i,j=1,\dots,d)$ is a $1$-cocycle whose
class in \v{C}ech cohomology $H^1(U_\bullet,\cO^\times_{\!X})$
corresponds to an invertible $\cO_{\!X}$-module $\cL$, with
a system of trivializations $\cO_{\!U_i}\isom\cL_{|U_i}$, for
$i=1,\dots,d$; especially, $\phi_i$ can be regarded as an
isomorphism $(\pi^{-1}\cL)_{|\pi^{-1}U_i}(r_i)\isom\cG_{|\pi^{-1}U_i}$
for every $i=1,\dots,d$. By inspecting the construction, it is
easily seen that these latter maps patch to a single isomorphism
of $\cO_{\P^n_{\!X}}$-modules $\pi^*\cL(r)\isom\cG$, as required.
\end{proof}

\begin{corollary}\label{cor_equivals}
Let $f:X\to S$ be a smooth morphism of finite presentation,
set\/ $\G_{m,X}:=X\times_S\Spec\,K^+[T,T^{-1}]$, and denote
by $\pi:\G_{m,X}\to X$ the natural morphism. Then the map
$$
\Pic\,X\to\Pic\,\G_{m,X}\qquad\cL\mapsto\pi^*\cL
$$
is an isomorphism.
\end{corollary}
\begin{proof} To start out, the map is injective, since
$\pi$ admits a section. Next, let $j:\G_{m,X}\to\P^1_{\!\!X}$
be the natural open immersion, and denote again by
$\pi:\P^1_{\!\!X}\to X$ the natural projection. According to
theorem \ref{th_cohereur} and proposition \ref{prop_hdim-finite}(i),
every invertible $\cO_{\G_{m,X}}$-module extends to an invertible
$\cO_{\P^1_{\!\!X}}$-module (namely the determinant of a coherent
extension). Since $j^*\cO_{\P^1_{\!\!X}}(n)\simeq\cO_{\G_{m,X}}$,
the claim follows from proposition \ref{prop_Mu}.
\end{proof}

\begin{remark} The proof of proposition \ref{prop_Mu} is based
on the argument given in \cite[Lecture 13, Prop.3]{Mu}. Of course,
statements related to corollary \ref{cor_equivals} abound in the
literature, and more general results are available: see {\em e.g.}
\cite{BaMu} and \cite{Sw}.
\end{remark}

\begin{proposition}\label{prop_equivals}
Let $f:X\to S$ be a smooth morphism of finite type, $j:U\to X$
an open immersion, and set $Z:=X\!\setminus\! U$.
Then:
\begin{enumerate}
\item
Every invertible $\cO_{\!U}$-module extends to an invertible
$\cO_{\!X}$-module.
\item
Suppose furthermore, that:
\begin{enumerate}
\item
For every point $y\in S$, the codimension of
$Z\cap f^{-1}(y)$ in $f^{-1}(y)$ is $\geq 1$, and
\item
The codimension of $Z\cap f^{-1}(\eta)$ in $f^{-1}(\eta)$
is $\geq 2$.
\end{enumerate}
Then the restriction functors:
\set\begin{equation}\label{eq_restrict}
\cO_{\!X}\bRflx\to\cO_{\!U}\bRflx \qquad \bPic\,X\to\bPic\,U
\end{equation}
are equivalences.
\end{enumerate}
\end{proposition}
\begin{proof} (i): Notice first that the underling space $|X|$ is
noetherian, since it is the union $|f^{-1}(s)|\cup|f^{-1}(\eta)|$
of two noetherian spaces. Especially, every open immersion is
quasi-compact, and $X$ is quasi-separated. We may then reduce easily
to the case where $U$ is dense in $X$, in which case the assertion
follows from propositions \ref{prop_extend-rflx}(i) and
\ref{prop_hdim-finite}(iv).

(ii): We begin with the following:

\begin{claim}\label{cl_faith-fct}
Under the assumptions of (ii), the functors \eqref{eq_restrict}
are faithful.
\end{claim}
\begin{pfclaim} Since $f$ is smooth, all the stalks of $\cO_{\!X}$
are reduced (\cite[Ch.IV, Prop.17.5.7]{EGA4}). Since every
point of $Z$ is specialization of a point of $U$, the claim
follows easily from remark \ref{rem_co-repres}.
\end{pfclaim}

Under assumptions (a) and (b), the conditions
of corollary \ref{cor_depth-cons} are fulfilled, and
one deduces that \eqref{eq_restrict} are full functors,
hence fully faithful, in view of claim \ref{cl_faith-fct}.
The essential surjectivity of the restriction functor for
reflexive $\cO_{\!X}$-modules is already known, by proposition
\ref{prop_extend-rflx}(i).
\end{proof}

\begin{lemma}\label{lem_rflx-pic}
Let $f:X\to S$ be a flat, finitely presented morphism,
and denote by $U\subset X$ the maximal open subset such that
the restriction $f_{|U}:U\to S$ is smooth. Suppose that the
following two conditions hold:
\begin{enumerate}
\alphaenu
\item
For every point $y\in S$, the fibre $f^{-1}(y)$ is geometrically
reduced.
\item
The generic fibre $f^{-1}(\eta)$ is geometrically normal.
\romanenu
\end{enumerate}
Then the restriction functor induces an equivalence of
categories (notation of \eqref{subsec_rk-trivias}) :
\set\begin{equation}\label{eq_restrict-gen-inv}
\bDiv\,X\isom\bPic\,U.
\end{equation}
\end{lemma}
\begin{proof} Let $\cF$ be any generically invertible reflexive
$\cO_{\!X}$-module; it follows from proposition
\ref{prop_hdim-finite}(iv) that $\cF_{|U}$ is an invertible
$\cO_{\!U}$-module, so the functor \eqref{eq_restrict-gen-inv}
is well defined.
Let $y\in S$ be any point; condition (a) says in particular that
the fibre $f^{-1}(y)$ is generically smooth over $\kappa(y)$
(\cite[Th.28.7]{Mat}). Then it follows from
\cite[Ch.IV, Th.17.5.1]{EGA4} that $U\cap f^{-1}(y)$ is
a dense open subset of $f^{-1}(y)$, and if
$x\in f^{-1}(y)\!\setminus\! U$, then the depth of
$\cO_{\!f^{-1}(y),x}$ is $\geq 1$, since the latter is a reduced
local ring of dimension $\geq 1$. Similarly, condition (b)
and Serre's criterion \cite[Ch.8, Th.23.8]{Mat}
say that for every $x\in f^{-1}(\eta)\!\setminus\! U$,
the local ring $\cO_{\!f^{-1}(\eta),x}$ has depth $\geq 2$.
Furthermore, \cite[Ch.IV, Prop.9.9.4]{EGAIV-3} implies that
the open immersion $j:U\to X$ is quasi-compact; summing up,
we see that all the conditions of corollary \ref{cor_depth-cons}
are fulfilled, so that $\cF=j_*j^*\cF$ for every reflexive
$\cO_{\!X}$-module $\cF$. This already means that
\eqref{eq_restrict-gen-inv} is fully faithful. To conclude,
it suffices to invoke proposition \ref{prop_extend-rflx}(i).
\end{proof}

\sset\subsubsection{}\label{subsec_semi-stable-curve}
Let $K(\sT)$ be the fraction field of the free polynomial
$K$-algebra $K[\sT]$. For every $\gamma\in\Gamma$ one can
define an extension of $|\cdot|$ to a {\em Gauss valuation\/}
$|\cdot|_{0,\gamma}:K(\sT)\to\Gamma$
(\cite[Ex.6.1.4(iii)]{Ga-Ra}). If $f(\sT):=\sum_{i=0}^da_i\sT^i$
is any polynomial, then
$|f(\sT)|_{0,\gamma}=\max\{|a_i|\cdot\gamma^i~|~i=0,\dots,d\}$.
We let $V(\gamma)$ be the valuation ring of $|\cdot|_{0,\gamma}$,
and set
$$
R(\gamma):=K[\sT,\sT^{-1}]\cap V(1)\cap V(\gamma).
$$
Since $K[\sT,\sT^{-1}]$ is a Dedekind domain, $R(\gamma)$ is
an intersection of valuation rings of the field $K(\sT)$,
hence it is a normal domain. By inspecting the definition
we see that $R(\gamma)$ consists of all the elements
of the form $f(\sT):=\sum^n_{i=-n}a_i\sT^i$, such that
$|a_i|\leq 1$ and $|a_i|\cdot\gamma^i\leq 1$ for every
$i=-n,\dots,n$. Suppose now that $\gamma\leq 1$, and
choose $c\in K^+$ with $|c|=\gamma$;
then every such $f(\sT)$ can be written uniquely in the
form $\sum^n_{i=0}a_i\sT^i+\sum^n_{j=1}b_j(c\sT^{-1})^j$,
where $a_i,b_j\in K^+$ for every $i,j\leq n$. Conversely,
every such expression yields an element of $R(\gamma)$.
In other words, we obtain a surjection of $K^+$-algebras
$K^+[\sX,\sY]\to R(\gamma)$ by the rule: $\sX\mapsto\sT$,
$\sY\mapsto c\sT^{-1}$. Obviously the kernel of this map
contains the ideal $(\sX\sY-c)$, and we leave to the reader
the verification that the induced map
\set\begin{equation}\label{eq_present-R}
K^+[\sX,\sY]/(\sX\sY-c)\to R(\gamma)
\end{equation}
is indeed an isomorphism.

\sset\subsubsection{}\label{subsec_gen-pt-over-s}
Throughout the following discussion, we shall assume
that $\gamma\leq 1$. Every $\delta\in\Gamma$ with
$\gamma\leq\delta\leq 1$, determines a prime ideal
$\fp(\delta):=\{f\in R(\gamma)~|~|f|_{0,\delta}<1\}\subset R(\gamma)$,
such that $\fm_KR(\gamma)\subset\fp(\delta)$.
Then it is easy to see that $R(\gamma)_{\fp(\delta)}\subset V(\delta)$,
and moreover :
$$
R(\gamma)_{\fp(\gamma)}=V(\gamma)\qquad R(\gamma)_{\fp(1)}=V(1)
$$
since $V(1)$ (resp. $V(\gamma)$) is already a localization
of $K^+[\sT]$ (resp. of $K^+[c\sT^{-1}]$).
In case $\gamma<1$, \eqref{eq_present-R} implies that
$R(\gamma)\otimes_{K^+}\kappa\simeq\kappa[\sX,\sY]/(\sX\sY)$, and
it follows easily that $\fp(1)$ and $\fp(\gamma)$ correspond
to the two minimal prime ideals of $\kappa[\sX,\sY]/(\sX\sY)$.
In case $\gamma=1$, we have
$R(1)\otimes_{K^+}\kappa\simeq\kappa[\sX,\sX^{-1}]$,
and again $\fp(1)$ corresponds to the generic point of
$\Spec\,\kappa[\sX,\sX^{-1}]$.  Notice that the natural
morphism
$$
f_\gamma:\T_{\!K}(\gamma):=\Spec\,R(\gamma)\to S
$$
restricts to a smooth morphism $f_\gamma^{-1}(\eta)\to\Spec\,K$;
moreover the closed fibre $f_\gamma^{-1}(s)$ is geometrically
reduced. Notice also that
$\T_{\!K}(\gamma)\times_S\Spec\,E^+\simeq\T_{\!E}(\gamma)$ for every
extension of valued fields $K\subset E$. In the following, we will
write just $\T(\gamma)$ in place of $\T_{\!K}(\gamma)$, unless we have
to deal with more than one base ring.

\begin{proposition}\label{prop_crit-reflex}
Keep the notation of \eqref{subsec_gen-pt-over-s}, and
let  $g:X\to\T(\gamma)$ be an {\'e}tale morphism, $\cF$ a
coherent $\cO_{\!X}$-module. Set $h:=f_\gamma\circ g:X\to S$
and denote by $i_s:h^{-1}(s)\to X$ the natural morphism.
Then $\cF$ is reflexive at the point $x\in h^{-1}(s)$
if and only if the following three conditions hold:
\begin{enumerate}
\alphaenu
\item
$\cF$ is $h$-flat at the point $x$.
\item
$\cF_{\!x}\otimes_{K^+}K$ is a reflexive $\cO_{\!X,x}\otimes_{K^+}K$-module.
\item
The $\cO_{\!h^{-1}(s),x}$-module $i^*_s\cF_{\!x}$ satisfies condition $S_1$
(see definition {\em\ref{def_Ass}(iii)}).
\romanenu
\end{enumerate}
\end{proposition}
\begin{proof} Suppose that $\cF$ is reflexive at the point $x$;
then it is easy to check that (a) and (b) hold. We prove (c): by
remark \ref{rem_co-repres} we can find a left exact sequence
$$
0\to\cF_{\!x}\stackrel{\alpha}{\to}\cO^{\oplus m}_{\!X,x}
\stackrel{\beta}{\to}\cO^{\oplus n}_{\!X,x}.
$$
Then $\Img\,\beta$ is a flat $K^+$-module, since it is a
submodule of the flat $K^+$-module $\cO^{\oplus n}_{\!X,x}$;
hence $\alpha\otimes_{K^+}\one_\kappa:
i^*_s\cF_{\!x}\to\cO^{\oplus m}_{\!h^{-1}(s),x}$
is still injective, so we are reduced to showing that
$h^{-1}(s)$ is a reduced scheme, which follows from
\cite[Ch.IV, Prop.17.5.7]{EGA4} and the fact that
$f_\gamma^{-1}(s)$ is reduced.

Conversely, suppose that conditions (a)--(c) hold.

\begin{claim}\label{cl_generic-vals}
Let $\xi$ be the generic point of an irreducible component
of $h^{-1}(s)$. Then:
\begin{enumerate}
\item
$\cF_{\!x}$ is a torsion-free $\cO_{\!X,x}$-module.
\item
$\cO_{\!X,\xi}$ is a valuation ring.
\item
Suppose that the closure of $\xi$ contains $x$. Then $\cF_\xi$
is a free $\cO_{\!X,\xi}$-module of finite rank.
\end{enumerate}
\end{claim}
\begin{pfclaim} (i): By (a), the natural map
$\cF_{\!x}\to\cF_{\!x}\otimes_{K^+}K$ is injective; since (b) implies
that $\cF_{\!x}\otimes_{K^+}K$ is a torsion-free $\cO_{\!X,x}$-module,
the same must then hold for $\cF_{\!x}$.

(ii) follows from proposition \ref{prop_integral-fp-ext}(ii).

(iii): Suppose that $x\in\overline{\{\xi\}}$. We derive easily
from (i) that $\cF_\xi$ is a torsion-free $\cO_{\!X,\xi}$-module,
so the assertion follows from (ii) and
\cite[Ch.VI, \S3, n.6, Lemma 1]{BouAC}.
\end{pfclaim}

By (b), the morphism
$\beta_{\cF\!,x}\otimes_{K^+}\one_K:\cF_{\!x}\otimes_{K^+}K\to
\cF^{\vee\vee}_{\!x}\otimes_{K^+}K$ is an isomorphism (notation
of \eqref{subsec_duals}); since $\cF$ is $h$-flat at $x$,
we deduce easily that $\beta_{\cF\!,x}$ is injective and
$C:=\Coker\,\beta_{\cF\!,x}$ is a torsion $K^+$-module.
To conclude, it remains only to show:

\begin{claim}\label{cl_C-is-flat} $C$ is a flat $K^+$-module.
\end{claim}
\begin{pfclaim}[] In view of lemma \ref{lem_Tor-crit}, it suffices
to show that $\Tor^{K^+}_1(C,\kappa(s))=0$. However, from the
foregoing we derive a left exact sequence
$$
\xymatrix{
0 \ar[r] & \Tor^{K^+}_1(C,\kappa(s)) \ar[r] &
\cF_{\!x}\otimes_{K^+}\kappa(s)
\ar[rr]^-{\beta_{\cF\!,x}\otimes_{K^+}\kappa(s)} & &
\cF^{\vee\vee}_{\!x}\otimes_{K^+}\kappa(s).
}
$$
We are thus reduced to showing that
$\beta_{\cF\!,x}\otimes_{K^+}\kappa(s)$ is an injective map.
In view of condition (c), it then suffices to prove that
$\beta_{\cF\!,\xi}$ is an isomorphism, whenever $\xi$ is the
generic point of an irreducible component of $h^{-1}(s)$
containing $x$. The latter assertion holds by virtue of
claim \ref{cl_generic-vals}(iii).
\end{pfclaim}
\end{proof}

\sset\subsubsection{}\label{subsec_basic-toys}
For a given $\rho\in\Gamma$, let us pick $a\in K\setminus\{0\}$
such that $|a|=\rho$; we define the fractional ideal
$I(\rho)\subset K[\sT,\sT^{-1}]$ as the $R(\gamma)$-submodule
generated by $\sT$ and $a$. The module $I(\rho)$ determines
a quasi-coherent $\cO_{\T(\gamma)}$-module $\cI(\rho)$.

\begin{lemma}\label{lem_goodrflx}
With the notation of \eqref{subsec_basic-toys}:
\begin{enumerate}
\item
$\cI(\rho)$ is a reflexive $\cO_{\T(\gamma)}$-module for
every $\rho\in\Gamma$.
\item
There exists a short exact sequence of $R(\gamma)$-modules:
$$
0\to I(\rho^{-1}\gamma)\to R(\gamma)^{\oplus 2}\to I(\rho)\to 0.
$$
\end{enumerate}
\end{lemma}
\begin{proof} To start with, let $a\in K\setminus\{0\}$
with $|a|=\rho$.

\begin{claim}\label{cl_symmetric}
If either $\rho\geq 1$ or $\rho\leq\gamma$, then $I(\rho)$
and $I(\rho^{-1}\gamma)$ are rank one, free $R(\gamma)$-modules.
\end{claim}
\begin{pfclaim} If $\rho\geq 1$ (resp. $\rho\leq\gamma$) then
$\rho^{-1}\gamma\leq\gamma$ (resp. $\rho^{-1}\gamma\geq 1$), hence
it suffices to show that $I(\rho)$ is free of rank one, in both
cases. Suppose first that $\rho\geq 1$; in this
case, multiplication by $a^{-1}$ yields an isomorphism
of $R(\gamma)$-modules $I(\rho)\isom R(\gamma)$.
Next, suppose that $\rho\leq\gamma$. Then $I(\rho)$
is the ideal $(\sT,a)$, where $a=c\cdot b$ for some $b\in K^+$
and $|c|=\gamma$. Therefore
$I(\rho)=(\sT,\sT\cdot(c\sT^{-1})\cdot b)=
\sT\cdot(1,c\sT^{-1}b)=\sT R(\gamma)$, and again $I(\rho)$ is
a free $R(\gamma)$-module of rank one.
\end{pfclaim}

In view of claim \ref{cl_symmetric}, we may assume that
$\gamma<\rho<1$. Let $(e_1,e_2)$ be the canonical basis
of the free $R(\gamma)$-module $R(\gamma)^{\oplus 2}$;
we consider the $R(\gamma)$-linear surjection
$\pi:R(\gamma)^{\oplus 2}\to I(\rho)$ determined
by the rule: $e_1\mapsto\sT$, $e_2\mapsto a$. Clearly
$\Ker\,\pi$ contains the submodule
$M(\rho)\subset R(\gamma)^{\oplus 2}$ generated by:
$$
f_1:=ae_1-\sT e_2 \quad \text{and} \quad
f_2:=c\sT^{-1}e_1-ca^{-1}e_2.
$$
Let $\cN$ be the quasi-coherent $\cO_{\T(\gamma)}$-module
associated with $N:=R(\gamma)^{\oplus 2}/M(\rho)$.

\begin{claim}\label{cl_N-is-flat} With the foregoing
notation:
\begin{enumerate}
\item
$N$ is a flat $K^+$-module.
\item
$K^+[c^{-1}]\otimes_{K^+}N$ is a free
$K^+[c^{-1}]\otimes_{K^+}R(\gamma)$-module of rank one.
\item
$K\otimes_{K^+}M(\rho)=K\otimes_{K^+}\Ker\,\pi$.
\end{enumerate}
\end{claim}
\begin{pfclaim}(i): We let $\gr_\bullet R(\gamma)$ be the
$\sT$-adic grading on $R(\gamma)$ ({\em i.e.}
$\gr_i R(\gamma)=\sT^i\cdot K\cap R(\gamma)$ for every $i\in\Z$),
and we define a compatible grading on $R(\gamma)^{\oplus 2}$ by
setting:
$\gr_i R(\gamma)^{\oplus 2}:=
(\gr_{i-1} R(\gamma)\cdot e_1)\oplus(\gr_i R(\gamma)\cdot e_2)$
for every $i\in\Z$. Since $f_1$ and $f_2$ are homogeneous elements,
we deduce by restriction a grading $\gr_\bullet M(\rho)$ on
$M(\rho)$, and a quotient grading $\gr_\bullet N$ on $N$,
whence a short exact sequence of graded $K^+$-modules:
$$
0\to\gr_\bullet M(\rho)\stackrel{\gr_\bullet j}{\longrightarrow}
\gr_\bullet R(\gamma)^{\oplus 2}\to\gr_\bullet N\to 0.
$$
However, by inspecting the definitions, it is easy to
see that $\gr_\bullet j$ is a split injective map of
free $K^+$-modules, hence $\gr_\bullet N$ is a free
$K^+$-module, and then the same holds for $N$.

(ii) is easy and shall be left to the reader.

(iii): Similarly, one checks easily that
$K\otimes_{K^+}I(\rho)$ is a free
$K\otimes_{K^+}R(\gamma)$-module of rank one;
then, by (ii) the quotient map
$K\otimes_{K^+}N\to K\otimes_{K^+}I(\rho)$ is necessarily
an isomorphism, whence the assertion.
\end{pfclaim}

\begin{claim}\label{cl_N-reflex}
$\cN$ is a reflexive $\cO_{\T(\gamma)}$-module.
\end{claim}
\begin{pfclaim} Since $\cN$ is coherent, it suffices to show
that $\cN_x$ is a reflexive $\cO_{\T(\gamma),x}$-module, for
every $x\in\T(\gamma)$ (lemma \ref{lem_refl-is-local}(ii)). Let
$y:=f_\gamma(x)$; we may then replace $\cN$ by its restriction
to $\T(\gamma)\times_SS(y)$, which allows to assume that $y=s$
is the closed point of $S$. In this case, we can apply the
criterion of proposition \ref{prop_crit-reflex} to the
morphism $f_\gamma:\T(\gamma)\to S$. We already know from
claim \ref{cl_N-is-flat}(i) that $\cN$ is $f_\gamma$-flat.
Moreover, by claim \ref{cl_N-is-flat}(ii) we see that
the restriction of $\cN$ to $f_\gamma^{-1}(\eta)$ is
reflexive. As $\gamma<\rho<1$, by inspecting the definition
and using the presentation \eqref{eq_present-R}, we deduce
an isomorphism :
$$
\kappa\otimes_{K^+}N\simeq\bar R{}^{\oplus 2}/(\sX e_2,\sY e_1)
\simeq(\bar R/\sX\bar R)\oplus(\bar R/\sY\bar R)
$$
where
$\bar R:=\kappa[\sX,\sY]/(\sX\sY)\simeq\kappa\otimes_{K^+}R(\gamma)$.
Thus, the $\bar R$-module $\kappa\otimes_{K^+}N$ satisfies condition
$S_1$, whence the claim.
\end{pfclaim}

It follows from claim \ref{cl_N-is-flat}(i,iii) that
the quotient map $N\to I(\rho)$ is an isomorphism, so
$\cI(\rho)$ is reflexive, by claim \ref{cl_N-reflex}.
Next, let us define an $R(\gamma)$-linear surjection
$\pi':R(\gamma)^{\oplus 2}\to M(\rho)$ by the rule:
$e_i\mapsto f_i$ for $i=1,2$. One checks easily that
$\Ker\,\pi'$ contains the submodule generated by the elements
$ca^{-1}e_1-Te_2$ and $ce_1-aTe_2$, and the latter is
none else than the module $M(\rho^{-1}\gamma)$, according
to our notation (notice that $\gamma<\rho^{-1}\gamma<1$).
We deduce a surjection of torsion-free
$R(\gamma)$-modules $I(\rho^{-1}\gamma)\isom
R(\gamma)^{\oplus 2}/M(\rho^{-1}\gamma)\to M(\rho)$, which
induces an isomorphism after tensoring by $K$, therefore
$I(\rho^{-1}\gamma)\isom M(\rho)$, which establishes (ii).
\end{proof}

\sset\subsubsection{}
Let $\T(\gamma)_\sm\subset\T(\gamma)$ be the largest open
subset which is smooth over $S$. Set $S_\gamma:=\Spec\,K^+/cK^+$;
it is easy to see that
$f_\gamma^{-1}(S\setminus S_\gamma)\subset\T(\gamma)_\sm$,
and for every $y\in S_\gamma$, the difference
$f_\gamma^{-1}(y)\setminus\T(\gamma)_\sm$ consists of a single
point.

\begin{proposition}\label{prop_calc-Pic-T}
Let $\Delta(\gamma)\subset\Gamma$ be the smallest convex
subgroup containing $\gamma$. Then there is a natural
isomorphism of groups:
$$
\Pic\,\T(\gamma)_\sm\isom\Delta(\gamma)/\gamma^\Z.
$$
\end{proposition}
\begin{proof} We consider the affine covering of
$\T(\gamma)_\sm$ consisting of the two open subsets
$$
U:=\Spec\,K^+[\sT,\sT^{-1}] \qquad \text{and} \qquad
V:=\Spec\,K^+[c\sT^{-1},c^{-1}\sT]
$$
with intersection $U\cap V=\Spec\,K^+[c^{-1},\sT,\sT^{-1}]$.
We notice that both $U$ and $V$ are $S$-isomorphic to $\G_{m,S}$,
and therefore
\set\begin{equation}\label{eq_pic-vanish}
\Pic\,U=\Pic\,V=0
\end{equation}
by corollary \ref{cor_equivals}. From \eqref{eq_pic-vanish},
a standard computation yields a natural isomorphism:
$$
\Pic\,\T(\gamma)_\sm\isom
\cO_{\T(\gamma)}(U)^\times\backslash
\cO_{\T(\gamma)}(U\cap V)^\times/\cO_{\T(\gamma)}(V)^\times.
$$
(Here, for a ring $A$, the notation $A^\times$ means the invertible
elements of $A$.)
However, $\cO_{\T(\gamma)}(U)^\times=(K^+)^\times\cdot(c\sT^{-1})^\Z$,
$\cO_{\T(\gamma)}(U\cap V)^\times=(K^+[c^{-1}])^\times\cdot\sT^\Z$
and $\cO_{\T(\gamma)}(V)=(K^+)^\times\cdot\sT^\Z$, whence the
contention.
\end{proof}

\sset\subsubsection{}\label{subsec_pic&reflex}
Proposition \ref{prop_calc-Pic-T} establishes a natural
bijection between the set of isomorphism classes of invertible
$\cO_{\T(\gamma)_\sm}$-modules and the set :
$$
]\gamma,1]:=\{\rho\in\Gamma~|~\gamma<\rho<1\}\cup\{1\}.
$$
On the other hand, lemma \ref{lem_rflx-pic} yields a
natural bijection between $\Pic\,\T(\gamma)_\sm$ and
the set of isomorphism classes of generically invertible
reflexive $\cO_{\T(\gamma)}$-modules. Furthermore, lemma
\ref{lem_goodrflx} provides already a collection of such
reflexive modules, and by inspection of the proof, we
see that the family of sheaves $\cI(\rho)$ is really
parametrized by the subset $]\gamma,1]$ (since the other
values of $\rho$ correspond to free
$\cO_{\T(\gamma)}$-modules of rank one). The two parametrizations
are essentially equivalent. Indeed, let $a\in K\setminus\{0\}$
be any element such that $\rho:=|a|\in]\gamma,1]$. With the
notation of the proof of proposition \ref{prop_calc-Pic-T},
we can define isomorphisms
$$
\phi:\cO_{\!U}\isom\cI(\rho)_{|U} \qquad \psi:\cO_V\isom\cI(\rho)_{|V}
$$
by letting: $\phi(1):=\sT$ and $\psi(1):=a$. To verify that
$\phi$ is an isomorphism, it suffices to remark that $\sT$
is a unit on $U$, so $\cI(\rho)_{|U}=\sT\cO_{\!U}=\cO_{\!U}$. Likewise,
on $V$ we can write $\sT=a\cdot(a^{-1}c)\cdot(c^{-1}\sT)$, so
$\cI(\rho)_{|V}=a\cO_V$, and $\psi$ is an isomorphism.
Hence $\cI(\rho)$ is isomorphic to the (unique)
$\cO_{\T(\gamma)}$-module whose global sections consist
of all the pairs $(s_U,s_V)\in\cO_{\!U}(U)\times\cO_V(V)$,
such that $\sT^{-1}s_{U|U\cap V}=a^{-1}s_{V|U\cap V}$.
Clearly, under the bijection of proposition
\ref{prop_calc-Pic-T}, the invertible sheaf
$\cI(\rho)_{|\T(\gamma)_\sm}$ corresponds to the class
of $\rho$ in $\Delta(\gamma)/\gamma^\Z$. In particular,
this shows that the reflexive $\cO_{\T(\gamma)}$-modules
$\cI(\rho)$ are pairwise non-isomorphic for
$\rho\in]\gamma,1]$, and that every reflexive
generically invertible $\cO_{\T(\gamma)}$-module is
isomorphic to one such $\cI(\rho)$.

\sset\subsubsection{}\label{subsec_also-local}
When $\gamma<1$ and $t\in\T(\gamma)$ is the singular point
of the closed fibre, the discussion of \eqref{subsec_pic&reflex}
also applies to describe the set $\bar\cohDiv(\cO_{\T(\gamma),t})$
of isomorphism classes of coherent reflexive fractional ideals
of $\cO_{\T(\gamma),t}$ (remark \ref{rem_classy-rflx}(ii)).
Indeed, any such module $M$ extends to a reflexive
$\cO_{\T(\gamma)}$-module (first one uses lemma
\ref{lem_rflx-on-lim}(ii.b) to extend $M$ to some quasi-compact
open subset $U\subset\T(\gamma)$, and then one may extend to the
whole of $\T(\gamma)$, via proposition \ref{prop_extend-rflx})(i).
Hence $M\simeq\cI(\rho)_t$ for some $\rho\in]\gamma,1]$. It follows
already that $\cohDiv(\cO_{\T(\gamma),t})$ is naturally an abelian
group, with multiplication law given by the rule :
$$
(M,N)\mapsto M\odot N:=j_*j^*(M\otimes_{\cO_{\T(\gamma),t}}N)
\qquad\text{for any two classes $M,N\in\Div(\cO_{\T(\gamma),t})$}
$$
where $j:\T(\gamma)_\sm\cap\Spec\,\cO_{\T(\gamma),t}\to
\Spec\,\cO_{\T(\gamma),t}$ is the natural open immersion. Indeed,
$M\odot N$ is reflexive (by proposition \ref{prop_extend-rflx}(i)
and corollary \ref{cor_depth-cons}), and the composition law $\odot$
is clearly associative and commutative, with $\cO_{\T(\gamma),t}$
as neutral element; moreover, for any $\rho\in\Delta(\gamma)$,
the class of $\cI(\rho)_t$ admits the inverse
$j_*((j^*\cI(\rho)_t)^\vee)$. Furthermore, the modules
$\cI(\rho)_t$ are pairwise non-isomorphic for $\rho\in]\gamma,1]$.
Indeed, using the group law $\odot$, the assertion follows once we
know that $\cI(\rho)_t$ is not trivial, whenever $\rho\in]\gamma,1[$.
However, from the presentation of lemma \ref{lem_goodrflx}(ii) one
sees that $\cI(\rho)_t\otimes_{\cO_{\T(\gamma),t}}\kappa(t)$ is a
two-dimensional $\kappa(t)$-vector space, so everything is clear.
Moreover, a simple inspection shows that the composition law
$\odot$ thus defined, agrees with the composition law of the
monoid $\cohDiv(\cO_{T(\gamma),t})$ given in remark
\ref{rem_classy-rflx}(ii). Summing, we get a natural group
isomorphism :
$$
\bar\cohDiv(\cO_{\T(\gamma),t})\isom\Delta(\gamma)/\gamma^{\Z}.
$$
The following theorem generalizes this classification to reflexive
modules of arbitrary generic rank.

\begin{theorem}\label{th_classify-rflx}
Let $g:X\to\T(\gamma)$ be an ind-\'etale morphism, $x\in X$
any point, $M$ a reflexive $\cO_{\!X,x}$-module.
Then there exist $\rho_1,\dots,\rho_n\in]\gamma,1]$ and
an isomorphism of $\cO_{\!X,x}$-modules:
$$
M\isom\bigoplus^n_{i=1} g^*\cI(\rho_i)_x.
$$
\end{theorem}
\begin{proof} Using lemma \ref{lem_rflx-on-lim}(ii.b), we are
easily reduced to the case where $g$ is \'etale.
Set $t:=g(x)$; first of all, if $t\in\T(\gamma)_\sm$,
then $X$ is smooth over $S$ at the point $t$, and consequently
$M$ (resp. $\cI(\rho)_t$) is a free $\cO_{\!X,x}$-module (resp.
$\cO_{\T(\gamma),t}$-module) of finite rank (by proposition
\ref{prop_hdim-finite}(iii)), so the assertion is obvious
in this case. Hence we may assume that $\gamma<1$ and $t$ is
the unique point in the closed fibre of
$\T(\gamma)\setminus\T(\gamma)_\sm$.
Next, let $K^{\sh+}$ be the strict henselization of $K^+$;
denote by $h:\T_{K^\sh}(\gamma)\to\T_K(\gamma)$ the natural
map, and set $X':=X\times_{\T_K(\gamma)}\T_{K^\sh}(\gamma)$.
Choose also a point $x'\in X'$ lying over $x$, and let
$t'\in\T_{K^\sh}(\gamma)$ be the image of $x'$; then $t'$ is
the unique point of $\T_{K^\sh}(\gamma)$ with $h(t')=t$.
We have a commutative diagram of
ring homomorphisms:
$$
\xymatrix{ \cO_{\T(\gamma),t} \ar[r]^-{g^\natural_x}
\ar[d]_{h^\natural_t} &
\cO_{\!X,x} \ar[d] \\
\cO_{\T_{K^\sh}(\gamma),t'} \ar[r] & \cO_{\!X',x'}
}
$$
whence an essentially commutative diagram of functors:
$$
\xymatrix{
\cO_{\T(\gamma),t}\bRflx \ar[r]^{g_x^*} \ar[d]_{h^*_t} &
\cO_{\!X,x}\bRflx \ar[d]^\beta \\
\cO_{\T_{K^\sh}(\gamma),t'}\bRflx \ar[r]^-\alpha &
\cO_{\!X',x'}\bRflx.
}
$$
Since $K$ and $K^\sh$ have the same value group, the discussion
in \eqref{subsec_also-local} shows that $h^*_t$ induces bijections
on the isomorphism classes of generically invertible modules.
On the other hand, proposition \ref{prop_descend-isom} implies
that the functor $\beta$ induces injections on isomorphism classes.
Consequently, in order to prove the theorem, we may replace the pair
$(X,x)$ by $(X',x')$, and assume from start that $K=K^\sh$.

In terms of the presentation \ref{eq_present-R} we can write
$\T(\gamma)\times_S\Spec\,\kappa(s)=Z_1\cup Z_2$, where $Z_1$
(resp. $Z_2$) is the reduced irreducible component on which $\sX$
(resp. $\sY$) vanishes. By inspecting the definitions, it is
easy to check that $Z_1\simeq\A^1_\kappa\simeq Z_2$ as
$\kappa$-schemes. Let $\xi_i$ be the generic point of $Z_i$,
for $i=1,2$; clearly $\{t\}=Z_1\cap Z_2$, hence
$W_i:=X\times_{\T(\gamma)}Z_i$ is non-empty and {\'e}tale over
$Z_i$, so $\cO_{W_i,x}$ is an integral domain, hence
$g^{-1}(\xi_i)$ contains exactly one point $\zeta_i$ that
specializes to $x$, for both $i=1,2$. To ease notation,
let us set $A:=\cO_{\!X,x}$. Then $\Spec\,A\otimes_{K^+}\kappa$
consists of exactly three points, namely $x$, $\zeta_1$ and
$\zeta_2$. Set $M(\zeta_i):=M\otimes_A\kappa(\zeta_i)$,
and notice that
$n:=\dim_{\kappa(\zeta_1)}M(\zeta_1)=\dim_{\kappa(\zeta_2)}M(\zeta_2)$,
since $M$ restricts to a locally free module over
$(\Spec\,A)_\sm$, the largest essentially smooth open $S$-subscheme
of $\Spec\,A$, which is connected. We choose a basis
$\bar e_1,\dots,\bar e_n$ (resp $\bar e'_1,\dots,\bar e'_n$)
for $M(\zeta_1)^\vee$ (resp. $M(\zeta_2)^\vee$), which we can then
lift to a system of sections
$e_1,\dots,e_n\in M^\vee_{\zeta_1}:=M^\vee\otimes_A\cO_{\!X,\zeta_1}$
(and likewise we construct a system $e'_1,\dots,e'_n\in M^\vee_{\zeta_2})$.
Let $\fp_i\subset A$ be the prime ideal corresponding to $\zeta_i$
($i=1,2$); after multiplication by an element of $A\setminus\fp_i$,
we may assume that $e_1,\dots,e_n\in M^\vee$ (and likewise for
$e'_1,\dots,e'_n$). Finally, we set $e''_i:=\sY e_i+\sX e'_i$ for
every $i=1,\dots,n$; it is clear that the system
$(e''_1,\dots,e''_n)$ induces bases of $M(\zeta_i)^\vee$ for both
$i=1,2$. We wish to consider the map:
$$
j:M\to A^{\oplus n}\qquad
m\mapsto(e''_1(m),\dots,e''_n(m)).
$$
Set $C:=\Coker\,j$, $I:=\Ann_A C$, and $B:=A/I$.

\begin{claim}\label{cl_about-j}
(i)\ \
The maps $M_{\zeta_i}\to\cO_{\!X,\zeta_i}^{\oplus n}$
induced by $j$ are isomorphisms.
\begin{enumerate}
\addenu
\item
$\Ker\,j=\Ker\,j\otimes_A\one_\kappa=0$.
\item
$B$ is a finitely presented $K^+$-module, and $C$
is a free $K^+$-module of finite rank.
\end{enumerate}
\end{claim}
\begin{pfclaim} (i): Using Nakayama's lemma, one deduces easily that
these maps are surjective; since $M_{\zeta_i}$ is a free
$\cO_{\!X,\zeta_i}$-module of rank $n$, they are also necessarily
injective.

(ii): Since $\cO_{\!X,x}$ is normal, $\{0\}$ is its only
associated prime; then the injectivity of $j$ (resp. of
$j\otimes_A\one_\kappa$) follows from (i), and the fact that
$M$ (resp. $M\otimes_A\kappa$) satisfies condition $S_1$, by
remark \ref{rem_co-repres} (resp. by proposition
\ref{prop_crit-reflex}).

(iii): First of all, since $A$ is coherent, $I$ is a finitely
generated ideal of $A$.
Let $f:\Spec\,B\to\Spec\,K^+$ be the natural morphism.
Since $C$ is a finitely presented $A$-module,
its support $Z$ is a closed subset of $\Spec\,A$.
From (i) we see that
$Z\cap\Spec\,A\otimes_{K^+}\kappa\subset\{x\}$; on the other hand,
$Z$ is also the support of the closed subscheme $\Spec\,B$
of $\Spec\,A$. Therefore
\set\begin{equation}\label{eq_Supp-B}
f^{-1}(s)\cap\Spec\,B\subset\{x\}.
\end{equation}
Since $K^+$ is henselian, it follows easily from \eqref{eq_Supp-B}
and \cite[Ch.IV, Th.18.5.11(c$^{\prime\prime}$)]{EGA4} that $B$ is
a finite $K^+$-algebra of finite presentation, hence also a finitely
presented $K^+$-module (claim \ref{cl_fin-present}). Since $C$ is a
finitely presented $B$-module, we conclude that $C$ is finitely
presented as $K^+$-module, as well. From (ii) we deduce a short
exact sequence: $0\to M\to A^{\oplus n}\to C\to 0$, and then
the long exact $\Tor$ sequence yields: $\Tor_1^A(C,\kappa)=0$
(cp. the proof of claim \ref{cl_C-is-flat}). We conclude
by \cite[Ch.II, \S3, n.2, Cor.2 of Prop.5]{BouAC}.
\end{pfclaim}

\begin{claim}\label{cl_B-closed-in-T}
The composed morphism: $\Spec\,B\to\Spec\,A\to\T(\gamma)$
is a closed immersion.
\end{claim}
\begin{pfclaim} Let $\fp\subset R(\gamma)$ be the maximal
ideal corresponding to $t$, so that
$\cO_{\T(\gamma),t}=R(\gamma)_\fp$. From claim
\ref{cl_about-j}(iii) we see that the natural
morphism $\psi:\,R(\gamma)_\fp\to B$ is finite.
Moreover, $A/\fp A\simeq R(\gamma)/\fp$, hence
$\psi\otimes_{R(\gamma)}\one_{R(\gamma)/\fp}$ is a surjection.
By Nakayama's lemma we deduce that $\psi$ is already a
surjection, {\em i.e.} the induced morphism
$\Spec\,B\to\Spec\,R(\gamma)_\fp$ is a closed immersion.
Let $J_\fp:=\Ker\,\psi$, $J:=R(\gamma)\cap J_\fp$
and $D:=R(\gamma)/J$. We are reduced to showing that
the induced map $D\to D_\fp$ is an isomorphism.
However, let $e\in R(\gamma)\setminus\fp$; since $D_\fp\simeq B$
is finite over $K^+$, we can find a monic polynomial
$P[T]\in K^+[T]$ such that $P(e^{-1})=0$, therefore an
identity of the form $1=e\cdot Q(e)$ holds in $D_\fp$
for some polynomial $Q(T)\in K^+[T]$. But then the same
identity holds already in the subring $D$, {\em i.e.}
the element $e$ is invertible in $D$, and the claim follows.
\end{pfclaim}

Now, $C$ is a finitely generated $B$-module,
hence also a finitely generated $R(\gamma)$-module, due
to claim \ref{cl_B-closed-in-T}. We construct a presentation
of $C$ in the following way. First of all, we have a
short exact sequence of
$R(\gamma)\otimes_{K^+}R(\gamma)$-modules:
$$
\underline{E}~:~0\to\Delta\to R(\gamma)\otimes_{K^+}R(\gamma)
\stackrel{\mu}{\to} R(\gamma)\to 0
$$
where $\mu$ is the multiplication map.
The homomorphism $R(\gamma)\to R(\gamma)\otimes_{K^+}R(\gamma)$
: $a\mapsto 1\otimes a$ fixes an $R(\gamma)$-module structure
on every $R(\gamma)\otimes_{K^+}R(\gamma)$-module (the {\em right\/}
$R(\gamma)$-module structure), and clearly $\underline{E}$ is
split exact, when regarded as a sequence of $R(\gamma)$-modules
via this homomorphism. Moreover, in terms of the presentation
\eqref{eq_present-R}, the $R(\gamma)\otimes_{K^+}R(\gamma)$-module
$\Delta$ is generated by the elements $\sX\otimes 1-1\otimes\sX$
and $\sY\otimes 1-1\otimes\sY$. Let $n$ be the rank of the free
$K^+$-module $C$ (claim \ref{cl_about-j}(iii)); there follows
an exact sequence
$$
\underline{E}\otimes_{R(\gamma)}C~:~
0\to\Delta\otimes_{R(\gamma)}C\to
R(\gamma)^{\oplus n}\to C\to 0
$$
which we may and do view as a short exact sequence of
$R(\gamma)$-modules, via the {\em left\/} $R(\gamma)$-module
structure induced by the restriction of scalars
$R(\gamma)\to R(\gamma)\otimes_{K^+}R(\gamma)$ :
$a\mapsto a\otimes 1$.
The elements $\sX,\sY\in R(\gamma)$ act as $K^+$-linear
endomorphisms on $C$; one can then find
bases $(b_i~|~i=1,\dots,n)$ and $(b'_i~|~i=1,\dots,n)$
of $C$, and elements $a_1,\dots,a_n\in K^+\setminus\{0\}$
such that $\sX b_i=a_ib'_i$ for every $i\leq n$. Since
$\sX\sY=c$ in $R(\gamma)$, it follows that $\sY b'_i=ca_i^{-1}b_i$
for every $i\leq n$. With this notation, it is clear
that $\Delta\otimes_{R(\gamma)}C$, with its left
$R(\gamma)$-module structure, is generated by the
elements:
$$
\sX\otimes b_i-1\otimes a_ib'_i\quad\text{and}\quad
\sY\otimes b'_i-1\otimes ca^{-1}_ib_i\quad (i=1,\dots,n).
$$
For every $i\leq n$, let $F_i$ be the $R(\gamma)$-module
generated freely by elements $(\eps_i,\eps'_i)$,
and $\Delta_i\subset F_i$ the submodule generated by
$\sX\eps_i-a_i\eps'_i$ and $\sY\eps'_i-ca^{-1}_i\eps_i$.
Moreover, let us write $b'_i=\sum^n_{j=1}u_{ij}b_j$
with unique $u_{ij}\in K^+$, let $F$ be the free
$R(\gamma)$-module with basis $(e_i~|~i=1,\dots,n)$,
and define $\phi:F\to\oplus^n_{i=1}F_i$ by the rule:
$e_i\mapsto\eps'_i-\sum^n_{j=1}u_{ij}\eps_j$ for every
$i\leq n$. We deduce a right exact sequence of
$R(\gamma)$-modules:
$$
F\oplus\bigoplus^n_{i=1}\Delta_i\stackrel{\psi_1}{\to}
\bigoplus^n_{i=1}F_i\stackrel{\psi_2}{\to} C\to 0
$$
where:
$$
\begin{aligned}
& \psi_1(f,d_1,\dots,d_n)=\phi(f)+(d_1,\dots,d_n) & &
\text{for every $f\in F$ and $d_i\in\Delta_i$} \\
& \psi_2(\eps_i)=b_i\quad\text{and}\quad\psi_2(\eps'_i)=b'_i & &
\text{for every $i=1,\dots,n$.}
\end{aligned}
$$
\begin{claim}\label{cl_psi_1-inj}
$\psi_1$ is injective.
\end{claim}
\begin{pfclaim} Let $L$ be the field of fractions of
$R(\gamma)$; since the domain of $\psi_1$ is a
torsion-free $K^+$-module, it suffices to verify that
$\psi\otimes_{R(\gamma)}\one_L$ is injective.
However, on the one hand $C\otimes_{R(\gamma)}L=0$,
and on the other hand, each $\Delta_i\otimes_{R(\gamma)}L$
is an $L$-vector space of dimension one, so the claim
follows by comparing dimensions.
\end{pfclaim}

By inspecting the definitions and the proof of lemma
\ref{lem_goodrflx}, one sees easily that
\set\begin{equation}\label{eq_Deltas}
\Delta_i\simeq I(|a^{-1}_ic|) \quad\text{for every $i\leq n$}.
\end{equation}
Moreover, by remark \eqref{rem_syzy}(iii), there
exist $p\in\N$ and an $R(\gamma)$-linear isomorphism:
\set\begin{equation}\label{eq_syz1}
\Ker\,\psi_2\isom R(\gamma)^{\oplus p}\oplus\Syz^1_{R(\gamma)}C.
\end{equation}
On the other hand, remark \eqref{rem_syzy}(iii) and claim
\ref{cl_about-j}(ii) also shows that there exist $q\in\N$
and an $A$-linear isomorphism:
\set\begin{equation}\label{eq_syz2}
M\isom A^{\oplus q}\oplus\Syz^1_AC.
\end{equation}
Combining \eqref{eq_syz1} and \eqref{eq_syz2} and using
lemma \ref{lem_syzy}, we deduce an $R(\gamma)^\he$-linear
isomorphism:
$$
\begin{aligned}
(R(\gamma)^\he)^{\oplus q}\oplus
(R(\gamma)^\he\otimes_{R(\gamma)}\Ker\,\psi_2) \isom\: &
(R(\gamma)^\he)^{\oplus p+q}\oplus
(R(\gamma)^\he\otimes_{R(\gamma)}\Syz^1_{R(\gamma)}C) \\
\isom\: & (R(\gamma)^\he)^{\oplus p+q}\oplus
\Syz^1_{R(\gamma)^\he}(R(\gamma)^\he\otimes_{R(\gamma)}C) \\
\isom\: & (R(\gamma)^\he)^{\oplus p+q}\oplus
(R(\gamma)^\he\otimes_A\Syz^1_AC) \\
\isom\; & (R(\gamma)^\he)^{\oplus p}\oplus(R(\gamma)^\he\otimes_AM).
\end{aligned}
$$
By claim \ref{cl_psi_1-inj} and \eqref{eq_Deltas} it follows
that $(R(\gamma)^\he)^{\oplus p}\oplus(R(\gamma)^\he\otimes_AM)$
is a direct sum of modules of the form
$R(\gamma)^\he\otimes_{R(\gamma)}I(\rho_i)$, for various
$\rho_i\in\Gamma$ (recall that $I(1)=R(\gamma)$). Notice
that every $I(\rho_i)$ is generically of rank one, hence
indecomposable. Then it follows from corollary \ref{cor_KRS-hens}
that $R(\gamma)^\he\otimes_AM$ is a direct sum of various
indecomposable $R(\gamma)^\he$-modules of the form
$R(\gamma)^\he\otimes_Ag^*\cI(\rho_i)_x$.
Finally, we apply proposition \ref{prop_descend-isom} to
conclude the proof of the theorem.
\end{proof}

The rest of this section shall be concerned with some
results that hold in the special case where the valuation
of $K$ has rank one.

\begin{theorem}\label{th_prim-dec-exist}
Suppose that $K$ is a valued field of rank one. Let $f:X\to S$ be
a finitely presented morphism, $\cF$ a coherent $\cO_{\!X}$-module.
Then :
\begin{enumerate}
\item
Every coherent proper submodule $\cG\subset\cF$ admits a
primary decomposition.
\item
$\Ass\,\cF$ is a finite set.
\end{enumerate}
\end{theorem}
\begin{proof} By lemma \ref{lem_a&b}, assertion (ii) follows from (i).
Using corollary \ref{cor_prim-dec-local}(ii) we reduce easily to the
case where $X$ is affine, say $X=\Spec\,A$ for a finitely presented
$K^+$-algebra $A$, and $\cF=M^\sim$, $\cG=N^\sim$ for some finitely
presented $A$-modules $N\subset M\neq 0$.
By considering the quotient $M/N$, we further reduce the proof to the
case where $N=0$. Let us choose a closed imbedding $i:X\to\Spec\,B$,
where $B:=K^+[T_1,\dots,T_r]$ is a free polynomial $K^+$-algebra;
by proposition \ref{prop_prim-dec&fin-map}, the submodule $0\subset\cF$
admits a primary decomposition if and only if the submodule
$0\subset i_*\cF$ does. Thus, we may replace $A$ by $B$, and
assume that $A=K^+[T_1,\dots,T_r]$, in which case we shall argue
by induction on $r$. Let $\xi$ denote the maximal point of
$f^{-1}(s)$; by claim \ref{cl_it-s-a-Gauss}, $V:=\cO_{\!X,\xi}$
is a valuation ring with value group $\Gamma$.
Then \cite[Lemma 6.1.14]{Ga-Ra} says that
$$
M_\xi\simeq V^{\oplus m}\oplus(V/b_1V)\oplus\cdots\oplus(V/b_kV)
$$
for some $m\in\N$ and elements $b_1,\dots,b_k\in\fm_K\!\setminus\!\{0\}$.
After clearing some denominators, we may then find a map
$\phi:M':=A^{\oplus m}\oplus(A/b_1A)\oplus\cdots\oplus(A/b_kA)\to M$
whose localization $\phi_\xi$ is an isomorphism.

\begin{claim}\label{cl_small-ass}
$\Ass\,A/bA=\{\xi\}$ whenever $b\in\fm_K\!\setminus\!\{0\}$.
\end{claim}
\begin{pfclaim} Clearly $\Ass\,A/bA\subset f^{-1}(s)$.
Let $x\in f^{-1}(s)$ be a non-maximal point. Thus, the
prime ideal $\fp\subset A$ corresponding to $x$ is not
contained in $\fm_KA$, {\em i.e.} there exists $a\in\fp$
such that $|a|_A=1$. Suppose by way of contradiction, that
$x\in\Ass\,A/bA$; then we may find $c\in A$ such that
$c\notin bA$ but $a^nc\in bA$ for some $n\in\N$. The
conditions translate respectively as the inequalities :
$$
|c|_A>|b|_A \qquad\text{and}\qquad |a^n|_A\cdot|c|_A=|a^nc|_A\leq|b|_A
$$
which are incompatible, since $|a^n|_A=1$.
\end{pfclaim}

Since $(\Ker\,\phi)_\xi=0$, claim \ref{cl_small-ass} and
proposition \ref{prop_misc-Ass}(ii) imply that
$\Ass\,\Ker\,\phi=\emptyset$; therefore $\Ker\,\phi=0$,
by lemma \ref{lem_Ass-Supp}(iii). Moreover, the submodule
$N_1:=A^{\oplus m}$ (resp. $N_2:=(A/b_1A)\oplus\cdots\oplus(A/b_kA)$)
of $M'$ is either $\xi$-primary (resp. $\{0\}$-primary), or else
equal to $M'$; in the first case, $0\subset M'$ admits the primary
decomposition $0=N_1\cap N_2$, and in the second case, either $0$
is primary, or else $M'=0$ has the empty primary decomposition.
Now, if $r=0$ then $M=M'$ and we are done. Suppose therefore that
$r>0$ and that the theorem is known for all integers $<r$. Set
$M'':=\Coker\,\phi$; clearly $\xi\notin\Supp\,M''$, so by
proposition \ref{prop_sh-exact_dec} we are reduced to showing
that the submodule $0\subset M''$ admits a primary decomposition.
We may then replace $M$ by $M''$ and assume from start that
$\xi\notin\Supp\,\cF$.
Let $\cF^t\subset\cF$ be the $K^+$-torsion submodule; clearly $\cF/\cF^t$
is a flat $K^+$-module, hence it is a finitely presented $\cO_{\!X}$-module
(proposition \ref{prop_Gruson-out}), so $\cF^t$ is a finitely generated
$\cO_{\!X}$-module, by \cite[Ch.10, \S.1, n.4, Prop.6]{BouAH}.
Hence we may find a non-zero $a\in\fm_K$ such that $a\cF^t=0$,
{\em i.e.} the natural map $\cF^t\to\cF/a\cF$ is injective.
Denote by $i:Z:=V(\Ann\,\cF/a\cF)\to X$ the closed immersion.

\begin{claim} There exists a finite morphism $g:Z\to\A^{r-1}_S$.
\end{claim}
\begin{pfclaim} Since $\xi\notin\Supp\,\cF$, we have $I\not\subset\fm_K A$,
{\em i.e.} we can find $b\in I$ such that $|b|_A=1$. Set $W:=V(a,b)\subset X$;
it suffices to exhibit a finite morphism $W\to\A^{r-1}_S$. By
\cite[\S5.2.4, Prop.2]{Bo-Gun}, we can find an automorphism
$\sigma:A/aA\to A/aA$ such that $\sigma(b)=u\cdot b'$,
where $u$ is a unit, and $b'\in(K^+/aK^+)[T_1,\dots,T_r]$
is of the form $b'=T_r^k+\sum_{j=0}^{k-1} a_jT^j_r$,
for some $a_0,\dots,a_{k-1}\in (K^+/aK^+)[T_1,\dots,T_{r-1}]$.
Hence the ring $A/(aA+bA)=A/(aA+b'A)$ is finite over
$K^+[T_1,\dots,T_{r-1}]$, and the claim follows.
\end{pfclaim}

\begin{claim}\label{cl_mod_by_a}
The $\cO_{\!X}$-submodule $0\subset\cF/a\cF$ admits a primary decomposition.
\end{claim}
\begin{pfclaim} Let $\cG:=(\cF/a\cF)_{|Z}$. By our inductive assumption
the submodule $0\subset g_*\cG$ on $\A^{r-1}_S$ admits a primary
decomposition; a first application of proposition \ref{prop_prim-dec&fin-map}
then shows that the submodule $0\subset\cG$ on $Z$ admits a primary
decomposition, and a second application proves the same for the submodule
$0\subset i_*\cG=\cF/a\cF$.
\end{pfclaim}

Finally, let $j:Y:=V(a)\subset X$ be the closed immersion, and set
$U:=X\setminus\! Y$. By construction, the natural map
$\underline\Gamma_Y\cF\to j_*j^*\cF$ is injective. Furthermore, $U$
is an affine noetherian scheme, so \cite[Th.6.8]{Mat} ensures that
the $\cO_{\!U}$-submodule $0\subset\cF_{|U}$ admits a primary decomposition.
The same holds for $j^*\cF$, in view of claim \ref{cl_mod_by_a}.
To conclude the proof, it remains only to invoke proposition
\ref{prop_prim_op&clo}.
\end{proof}

\begin{corollary}\label{cor_prim-dec-exist}
Let $K$ be a valued field of rank one, $f:X\to S$ a finitely
presented morphism, $\cF$ a coherent $\cO_{\!X}$-module, and
$\underline\cI:=(\cI_\lambda~|~\lambda\in\Lambda)$ a cofiltered
family such that :
\begin{enumerate}
\alphaenu
\item
$\cI_\lambda\subset\cO_{\!X}$ is a coherent ideal for every
$\lambda\in\Lambda$.
\item
$\cI_\lambda\cdot\cI_\mu\in\underline\cI$ whenever
$\cI_\lambda,\cI_\mu\in\underline\cI$.
\romanenu
\end{enumerate}
Then the following holds :
\begin{enumerate}
\item
$\cF_\lambda:=\Ann_\cF(\cI_\lambda)\subset\cF$ is a submodule of
finite type for every $\lambda\in\Lambda$.
\item
There exists $\lambda\in\Lambda$ such that $\cF_\mu\subset\cF_\lambda$
for every $\mu\in\Lambda$.
\end{enumerate}
\end{corollary}
\begin{proof} We easily reduce to the case where $X$ is affine,
say $X=\Spec\,A$ with $A$ finitely presented over $K^+$, and for
each $\lambda\in\Lambda$ the ideal
$I_\lambda:=\Gamma(X,\cI_\lambda)\subset A$ is finitely generated.

(i): For given $\lambda\in\Lambda$, let $f_1,\dots,f_k$ be a
finite system of generators of $I_\lambda$; then $\cF_\lambda$
is the kernel of the map $\phi:\cF\to\cF^{\oplus k}$ defined by
the rule : $m\mapsto(f_1m,\dots,f_km)$, which is finitely
generated because $A$ is coherent.

(ii): By theorem \ref{th_prim-dec-exist} we can find primary
submodules $\cG_1,\dots,\cG_n\subset\cF$ such that
$\cG_1\cap\cdots\cap\cG_n=0$; for every $i\leq n$ and
$\lambda\in\Lambda$ set $\cH_i:=\cF/\cG_i$ and
$\cH_{i,\lambda}:=\Ann_{\cH_i}(\cI_\lambda)$.
Since the natural map $\cF\to\oplus^n_{i=1}\cH_i$
is injective, we have :
$$
\cF_\lambda=\cF\cap(\cH_{1,\lambda}\oplus\cdots\oplus\cH_{n,\lambda})
\qquad\text{for every $\lambda\in\Lambda$}.
$$
It suffices therefore to prove that, for every $i\leq n$,
there exists $\lambda\in\Lambda$ such that
$\cH_{i,\mu}=\cH_{i,\lambda}$ for every $\mu\in\Lambda$.
Say that $\cH_i$ is $\fp$-primary, for some prime ideal
$\fp\subset A$; suppose now that there exists $\lambda$
such that $I_\lambda\subset\fp$; since $I_\lambda$ is
finitely generated, we deduce that $I_\lambda^n\cH_i=0$
for $n\in\N$ large enough; from (b) we see that
$\cI^n_\lambda\in\underline\cI$, so (ii) holds in this case.
In case $I_\lambda\not\subset\fp$ for every $\lambda\in\Lambda$,
we have $\Ann_{\cH_i}(\cI_\lambda)=0$ for every
$\lambda\in\Lambda$, so (ii) holds in this case as well.
\end{proof}

\begin{corollary}\label{cor_integral-fp-ext}
Suppose $A$ is a valuation ring with value group
$\Gamma_{\!\!A}$, and $\phi:K^+\to A$ is an essentially
finitely presented local homomorphism from a valuation ring
$K^+$ of rank one. Then:
\begin{enumerate}
\item
If\/ the valuation of $K$ is not discrete, $\phi$ induces
an isomorphism $\Gamma\isom\Gamma_{\!\!A}$.
\item
If\/ $\Gamma\simeq\Z$ and $\phi$ is flat, $\phi$ induces an
inclusion $\Gamma\subset\Gamma_{\!\!A}$, and $(\Gamma_{\!\!A}:\Gamma)$
is finite.
\end{enumerate}
\end{corollary}
\begin{proof} Suppose first that $\Gamma\simeq\Z$; then $A$ is
noetherian, hence $\Gamma_{\!\!A}$ is discrete of rank one as well,
and the assertion follows easily. In case $\Gamma$ is not discrete,
we claim that $A$ has rank $\leq 1$. Indeed, suppose by way of
contradiction, that the rank of $A$ is higher than one, and let
$\fm_A\subset A$ be the maximal ideal; then we can find
$a,b\in\fm_A\!\setminus\!\{0\}$ such that $a^{-i}b\in A$ for
every $i\in\N$. Let us consider the $A$-module $M:=A/bA$;
we notice that $\Ann_M(a^i)=a^{-i}bA$ form a strictly
increasing sequence of ideals, contradicting corollary
\ref{cor_prim-dec-exist}(ii). Next we claim that $\phi$
is flat. Indeed, suppose this is not the case; then
$A$ is a $\kappa$-algebra. Let $\fp\subset A$ be the maximal
ideal; by lemma \ref{lem_integral-fp-ext}, we may assume that
$A/\fp$ is a finite extension of $\kappa$.
Now, choose any finitely presented quotient $\bar A$
of $A$ supported at $\fp$; it follows easily that $\bar A$
is integral over $K^+$, hence it is a finitely presented
$K^+$-module by proposition \ref{prop_integral-fp-ext}(i),
and its annihilator contains $\fm_K$, which is absurd in view
of \cite[Lemma 6.1.14]{Ga-Ra}.
Next, since $\Gamma$ is not discrete, one sees easily
that $\fm_K A$ is a prime ideal, and in light of the foregoing,
it must then be the maximal ideal. Now the assertion follows
from proposition \ref{prop_integral-fp-ext}(ii).
\end{proof}

\sset\subsubsection{}\label{subsec_back-to-Rees}
Let $K$ be a valued field of rank one, $B$ a finitely
presented $K^+$-algebra, and consider a pair
$\underline A:=(A,\Fil_\bullet A)$ consisting of a
$B$-algebra and a $B$-algebra filtration on $A$; let
also $M$ be a finitely generated $A$-module, and
$\Fil_\bullet M$ a good $\underline A$-filtration on $M$
(see definition \ref{def_good-filtr}). By definition,
this means that the Rees module $\sR(\underline M)_\bullet$
of the filtered $\underline A$-module
$\underline M:=(M,\Fil_\bullet M)$ is finitely generated over
the graded Rees $B$-algebra $\sR(\underline A)_\bullet$.
We have :

\begin{proposition}\label{prop_Rees}
In the situation of \eqref{subsec_back-to-Rees}, suppose
that $A$ is a finitely presented $B$-algebra, $M$ is a
finitely presented $A$-module, and $\Fil_\bullet A$ is a
good filtration. Then :
\begin{enumerate}
\item
$\sR(\underline A)_\bullet$ is a finitely presented $B$-algebra,
and $\sR(\underline M)_\bullet$ is a finitely presented
$\sR(\underline A)_\bullet$-module.
\item
If furthermore, $\Fil_\bullet A$ is a positive filtration
(see definition {\em\ref{def_good-filtr}(ii)}), then $\Fil_iM$
is a finitely presented $B$-module, for every $i\in\Z$.
\end{enumerate}
\end{proposition}
\begin{proof} By lemma \ref{lem_good-filtr}, there exists a
finite system of generators $\bm:=(m_1,\dots,m_n)$ of $M$,
and a sequence of integers $\bk:=(k_1,\dots,k_n)$ such that
$\Fil_\bullet M$ is of the form \eqref{eq_standard-filtr}.

(i): Consider first the case where $A$ is a free $B$-algebra of
finite type, say $A=B[t_1,\dots,t_p]$, such that $\Fil_\bullet A$
is the good filtration associated with the system of generators
$\bt:=(t_1,\dots,t_p)$ and the sequence of integers
$\br:=(r_1,\dots,r_p)$, and moreover $M$ is a free $A$-module
with basis $\bm$. Then $\sR(\underline A)_\bullet$
is also a free $B$-algebra of finite type (see example
\ref{ex_Rees-free}). Moreover, for every $j\leq n$, let
$M_j\subset M$ be the $A$-submodule generated by $m_j$,
and denote by $\Fil_\bullet M_j$ the good $\underline A$-filtration
associated with the pair $(\{m_j\},\{k_j\})$; clearly
$\Fil_\bullet M=\Fil_\bullet M_1\oplus\cdots\oplus\Fil_\bullet M_n$,
therefore $\sR(\underline M)_\bullet=
\sR(\underline M_1)_\bullet\oplus\cdots\oplus\sR(\underline M_n)_\bullet$,
where $\underline M_j:=(M_j,\Fil_\bullet M_j)$ for every $j\leq n$.
Obviously each
$\sR(\underline A)_\bullet$-module $\sR(\underline M_j)_\bullet$
is free of rank one, so the proposition follows in this case.

Next, suppose that $A$ is a free $B$-algebra of finite type,
and $M$ is arbitrary. Let $F$ be a free $A$-module of rank $n$,
$\be:=(e_1,\dots,e_n)$ a basis of $F$, and define an $A$-linear
surjection $\phi:F\to M$ by the rule $e_j\mapsto m_j$ for every
$j\leq n$. Then $\phi$ is even a map of filtered $\underline A$-modules,
provided we endow $F$ with the good $\underline A$-filtration
$\Fil_\bullet F$ associated with the pair $(\be,\bk)$. More
precisely, let $N:=\Ker\,\phi$; then the filtration $\Fil_\bullet M$
is induced by $\Fil_\bullet F$, meaning that $\Fil_iM:=(N+\Fil_iF)/N$
for every $i\in\Z$. Obviously the natural map
$\sR(\underline F)_\bullet\to\sR(\underline M)_\bullet$
is surjective, and its kernel is the Rees module
$\sR(\underline N)_\bullet$ corresponding to the filtered
$\underline A$-module $\underline N:=(N,\Fil_\bullet N)$
with $\Fil_iN:=N\cap\Fil_iF$ for every $i\in\Z$. Let
$\bx:=(x_1,\dots,x_s)$ be a finite system of generators of $N$,
and choose a sequence of integers $\bj:=(j_1,\dots,j_s)$ such
that $x_i\in\Fil_{j_i}N$ for every $i\leq s$; denote by
$\underline L$ the filtered $\underline A$-module associated
as in \eqref{subsec_good-filtr}, with the $A$-module $N$ and the
pair $(\bx,\bj)$. Thus, $\sR(\underline L)_\bullet$ is a finitely
generated graded $\sR(\underline A)_\bullet$-submodule of
$\sR(\underline N)_\bullet$. To ease notation, set
$\bar N:=\sR(\underline N)_\bullet/\sR(\underline L)_\bullet$
and $\bar F:=\sR(\underline F)_\bullet/\sR(\underline L)_\bullet$.
Recall (definition \ref{def_Rees}(iii)) that $\sR(\underline A)_\bullet$
is a graded $B$-subalgebra of $A[U,U^{-1}]$; then we have :

\begin{claim}\label{cl_annul-Rees}
$\bar N=\bigcup_{n\in\N}\Ann_{\bar F}(U^n)$.
\end{claim}
\begin{pfclaim} It suffices to consider a homogeneous element
$y:=U^iz\in\sR(\underline F)_i$, for some some $z\in\Fil_iF$.
Suppose that $U^ny\in\sR(\underline L)_{i+n}$;
especially, $U^{i+n}z\in\sR(\underline N)_{i+n}$, so $z\in N$,
hence $y\in\sR(\underline N)_i$. Conversely, suppose that
$y\in\sR(\underline N)_i$; write $z=x_1a_1+\cdots+x_sa_s$
for some $a_1,\dots,a_s\in A$. Say that $a_r\in\Fil_{b_r}A$
for every $r\leq s$, and set $l:=\max(i,b_1+j_1,\dots,b_s+j_s)$.
Then $U^{l-i}y=U^lz\in\sR(\underline L)_l$.
\end{pfclaim}

By the foregoing, $\sR(\underline F)_\bullet$ is finitely
presented over $\sR(\underline A)_\bullet$. It follows that
$\bar N$ is finitely generated (corollary \ref{cor_prim-dec-exist}),
hence $\sR(\underline N)_\bullet$ is a finitely generated
$\sR(\underline A)_\bullet$-module, so finally
$\sR(\underline M)_\bullet$ is finitely presented over
$\sR(\underline A)_\bullet$.

Lastly, consider the case of an arbitrary finitely presented
$B$-algebra $A$, endowed with the good $B$-algebra filtration
associated with a system of generators $\bx:=(x_1,\dots,x_p)$
and the sequence of integers $\br$. We map the free $B$-algebra
$P:=B[t_1,\dots,t_p]$ onto $A$, by the rule: $t_j\mapsto x_j$
for every $j\leq p$. This is even a map of filtered algebras,
provided we endow $P$ with the good filtration $\Fil_\bullet P$
associated with the pair $(\bt,\br)$; more precisely, $\Fil_\bullet P$
induces the filtration $\Fil_\bullet A$ on $A$, hence $\Fil_\bullet A$
is a good $\underline P$-filtration. By the foregoing
case, we then deduce that $\sR(\underline A)_\bullet$ is a
finitely presented $\sR(\underline P)_\bullet$-module (where
$\underline P:=(P,\Fil_\bullet P)$), hence also a finitely presented
$K^+$-algebra. Moreover, $M$ is a finitely presented $P$-module,
and clearly $\Fil_\bullet M$ is good when regarded as a
$\underline P$-filtration, hence $\sR(\underline M)_\bullet$
is a finitely presented $\sR(\underline P)_\bullet$-module,
so also a finitely presented $\sR(\underline A)_\bullet$-module.

(ii): The positivity condition implies that $\sR(\underline A)_0=B$;
then, taking into account (i), the assertion is just a special case
of proposition \ref{prop_four-year-later}(iii).
\end{proof}

\begin{theorem}[Artin-Rees lemma]\label{th_Rees}
Let $K$ be a valued field of rank one, $A$ an essentially
finitely presented $K^+$-algebra, $I\subset A$ an ideal of
finite type, $M$ a finitely presented $A$-module,
$N\subset M$ a finitely generated submodule.
Then there exists an integer $c\in\N$ such that :
$$
I^nM\cap N=I^{n-c}(I^cM\cap N)
\qquad\text{for every $n\geq c$}.
$$
\end{theorem}
\begin{proof} We reduce easily to the case where $A$ is a
finitely presented $K^+$-algebra; then, with all the work done
so far, we only have to repeat the argument familiar from the
classical noetherian case. Indeed, let us define a filtered
$K^+$-algebra $\underline A:=(A,\Fil_\bullet A)$ by the rule :
$\Fil_nA:=I^{-n}$ if $n\leq 0$, and $\Fil_nA:=A$ otherwise;
also let $\underline M:=(M,\Fil_\bullet M)$ be the filtered
$\underline A$-module such that $\Fil_iM:=M\cdot\Fil_iA$ for
every $i\in\Z$. We endow $Q:=M/N$ with the filtration
$\Fil_\bullet Q$ induced from $\Fil_\bullet M$,
so that the natural projection $M\to M/N$ yields a map
of filtered $\underline A$-modules
$\underline M\to\underline Q:=(Q,\Fil_\bullet Q)$.
By proposition \ref{prop_Rees}, there follows a surjective map
of finitely presented graded $\sR(\underline A)_\bullet$-modules :
$\pi_\bullet:\sR(\underline M)_\bullet\to\sR(\underline Q)_\bullet$,
whose kernel in degree $k\leq 0$ is the $A$-module $I^{-k}M\cap N$.
Hence we can find a finite (non-empty) system $f_1,\dots,f_r$
of generators for the $\sR(\underline A)_\bullet$-module
$\Ker\,\pi_\bullet$; clearly we may assume that each $f_i$
is homogeneous, say of degree $d_i$; we may also suppose
that $d_i\leq 0$ for every $i\leq r$, since $\Ker\,\pi_0$
generates $\Ker\,\pi_n$, for every $n\geq 0$.
Let $d:=\min\,(d_i~|~i=1,\dots,r)$; by inspecting the
definitions, one verifies easily that
$$
\sR(\underline A)_i\cdot\Ker\,\pi_j\subset
\sR(\underline A)_{i+j-d}\cdot\Ker\,\pi_d
\qquad\text{whenever $i+j\leq d$ and $0\geq j\geq d$}.
$$
Therefore $\Ker\,\pi_{d+k}=\sR(\underline A)_k\cdot\Ker\,\pi_d$
for every $k\leq 0$, so the assertion holds with $c:=-d$.
\end{proof}

\sset\subsubsection{}\label{subsec_Rees-complt}
In the situation of theorem \ref{th_Rees}, let $M$
be any finitely presented $A$-module; we denote by $M^\wedge$
the $I$-adic completion of $M$, which is an $A^\wedge$-module.

\begin{corollary}\label{cor_cption-exact}
With the notation of \eqref{subsec_Rees-complt},
the following holds :
\begin{enumerate}
\item
The functor $M\mapsto M^\wedge$ is exact on the
category of finitely presented $A$-modules.
\item
For every finitely presented $A$-module $M$, the natural
map $M\otimes_AA^\wedge\to M^\wedge$ is an isomorphism.
\item
$A^\wedge$ is flat over $A$. If additionally $I\subset\rad\,A$
(the Jacobson radical of $A$), then $A^\wedge$ is faithfully
flat over $A$.
\end{enumerate}
\end{corollary}
\begin{proof} (i): Let $N\subset M$ be any injection of
finitely presented $A$-modules; by theorem \ref{th_Rees},
the $I$-adic topology on $M$ induces the $I$-adic topology
on $N$. Then the assertion follows from proposition
\ref{prop_replaces-Mat-Th.8.1}(i,v).

(ii): Choose a presentation
$A^{\oplus p}\to A^{\oplus q}\to M\to 0$. By (i) we deduce
a commutative diagram with exact rows :
$$
\xymatrix{
A^{\oplus p}\otimes_AA^\wedge \ar[r] \ar[d] &
A^{\oplus q}\otimes_AA^\wedge \ar[r] \ar[d] &
M\otimes_AA^\wedge \ar[r] \ar[d] & 0 \\
(A^{\oplus p})^\wedge \ar[r] & (A^{\oplus q})^\wedge \ar[r] &
M^\wedge \ar[r] & 0
}$$
and clearly the two left-most vertical arrows are isomorphisms.
The claim follows.

(iii): The first assertion means that the functor
$M\mapsto M\otimes_AA^\wedge$ is exact; this follows from
(i) and (ii), via a standard reduction to the case where $M$
is finitely presented. Suppose next that $I\subset\rad\,A$;
to conclude, it suffices to show that the image of the natural
map $\Spec\,A^\wedge\to\Spec\,A$  contains the maximal spectrum
$\Max\,A$ (\cite[Th.7.3]{Mat}). Since the natural map
$\Max\,A/I\to\Max\,A$ is a bijection, the latter assertion
follows from the :

\begin{claim}\label{cl_was-skipped}
For every finitely presented $A$-module $M$, the natural map
$i_M:M\to M^\wedge$ induces an isomorphism
$M/IM\isom M^\wedge/IM^\wedge$.
\end{claim}
\begin{pfclaim}[] From (i) we get a natural identification :
$M^\wedge/(IM)^\wedge\isom(M/IM)^\wedge\isom M/IM$, whose inverse
is the map induced by $i_M$. However, by (ii) the image of
$(IM)^\wedge$ in $M^\wedge$ is the same as the image of
$IM\otimes_AA^\wedge$, which is the same as the image of
$IM^\wedge$.
\end{pfclaim}
\end{proof}

\begin{theorem}\label{th_omega-coh}
Let $K$ be a valued field of rank one, $B\to A$ a map of
finitely presented $K^+$-algebras. Then :
\begin{enumerate}
\item
If $M$ is an $\omega$-coherent $A$-module, $M$ is
$\omega$-coherent as a $B$-module.
\item
If $J$ is a coh-injective $B$-module, and $I\subset B$ is a
finitely generated ideal, we have :
\begin{enumerate}
\item
$\Hom_B(A,J)$ is a coh-injective $A$-module.
\item
$\bigcup_{n\in\N}\Ann_J(I^n)$ is a coh-injective $B$-module.
\end{enumerate}
\item
If $(J_n,\phi_n~|~n\in\N)$ is a direct system consisting of
coh-injective $B/I^n$-modules $J_n$ and $B$-linear maps
$\phi_n:J_n\to J_{n+1}$ (for every $n\in\N$), then
$\colim_{n\in\N}J_n$ is a coh-injective $B$-module.
\end{enumerate}
\end{theorem}
\begin{proof}(See \eqref{subsec_coh-inject} for the generalities
on coh-injective and $\omega$-coherent modules.)

(i): We reduce easily to the case where $M$ is
finitely presented over $A$. Let $x_1,\dots,x_p$ be a
system of generators for the $B$-algebra $A$, and
$m_1,\dots,m_n$ a system of generators for the $A$-module
$M$. We let $\underline A:=(A,\Fil_\bullet A)$, where
$\Fil_\bullet A$ is the good $B$-algebra filtration associated
with the pair $\bx:=(x_1,\dots,x_p)$ and $\br:=(1,\dots,1)$;
likewise, let $\underline M$ be the filtered $\underline A$-module
defined by the good $\underline A$-filtration on $M$
associated with $\bm:=(m_1,\dots,m_n)$ and $\bk:=(0,\dots,0)$
(see definition \ref{def_good-filtr}). Then claim follows easily
after applying proposition \ref{prop_Rees}(ii) to the filtered
$B$-algebra $\underline A$ and the filtered $\underline A$-module
$\underline M$.

(ii.a): Let $N\subset M$ be two coherent $A$-modules,
and $\phi:N\to\Hom_B(A,J)$ an $A$-linear map. According to claim
\ref{cl_right-to-forget}, $\phi$ corresponds by adjunction to a unique
$B$-linear map $\bar\phi:N\to J$; on the other hand by (i), $M$,
and $M/N$ are $\omega$-coherent $B$-modules, hence $\bar\phi$
extends to a $B$-linear map $\bar\psi:M\to J$ (lemma
\ref{lem_omega-coh}). Under the adjunction, $\bar\psi$ corresponds
to an $A$-linear extension $\psi:M\to\Hom_B(A,J)$ of $\phi$.

(iii): Let $M\subset N$ be two finitely presented $B$-modules,
$f:M\to J:=\colim_{n\in\N}J_n$ a $B$-linear map. Since $M$ is
finitely generated, $f$ factors through a map $f_n:M\to J_n$
and the natural map $J_n\to J$, provided $n$ is large enough
(\cite[Prop.2.3.16(ii)]{Ga-Ra}).
By theorem \ref{th_Rees} there exists $c\in\N$ such that
$I^{n+c}N\cap M\subset I^nM$. Hence
$$
f_{n+c}:=\phi_{n+c-1}\circ\cdots\circ\phi_n\circ f_n:M\to J_{n+c}
$$
factors through a unique $B/I^{n+c}$-linear map
$\bar f_{n+c}:\bar M:=M/(I^{n+c}N\cap M)\to J_{n+c}$; since $I^{n+c}$
is finitely generated, $\bar M$ is a coherent submodule of the
coherent $B$-module $\bar N:=N/I^{n+c}N$, therefore $\bar f_{n+c}$
extends to a map $\bar g:\bar N\to J_{n+c}$. The composition of
$\bar g:\bar N\to J_{n+c}$ with the projection $N\to\bar N$ and
the natural map $J_{n+c}\to J$, is the sought extension of $f$.

(ii.b): Letting $A:=B/I^n$ in (ii.a), we deduce that
$J_n:=\Ann_J(I^n)$ is a coh-injective $B/I^n$-module, for
every $n\in\N$. Then the assertion follows from (iii).
\end{proof}

\subsection{Local duality}\label{sec_loc-duality}
Throughout this section we let $(K,|\cdot|)$ be a valued
field of rank one. We shall continue to use the general
notation of \eqref{sec_sch-val-rings}.

\sset\subsubsection{}\label{subsec_ext-and-Gamma_Z}
Let $A$ be finitely presented $K^+$-algebra,
$I\subset A$ an ideal generated by a finite system
$\bff:=(f_i~|~i=1,\dots,r)$, and denote by
$i:Z:=V(I)\to X:=\Spec\,A$ the natural closed immersion.
Let also $\cI\subset\cO_{\!X}$ be the ideal arising from $I$.
For every $n\geq 0$ there is a natural epimorphism
$$
i_*i^{-1}\cO_{\!X}\to\cO_{\!X}/\cI^n
$$
whence natural morphisms in $\sD(A\Mod)$ :
\set\begin{equation}\label{eq_ext-and-gamma_Z}
R\Hom^\bullet_{\cO_{\!X}}(\cO_{\!X}/\cI^n,\cF^\bullet)\to
R\Hom^\bullet_{\cO_{\!X}}(i_*i^{-1}\cO_{\!X},\cF^\bullet)\isom
R\Gamma_{\!\!Z}\cF^\bullet
\end{equation}
for any bounded below complex $\cF^\bullet$ of $\cO_{\!X}$-modules.

\begin{theorem}\label{th_calculate_Z}
In the situation of \eqref{subsec_ext-and-Gamma_Z}, let
$M^\bullet$ be any object of $\sD^+(A\Mod)$.
Then \eqref{eq_ext-and-gamma_Z} induces natural isomorphisms :
$$
\colim_{n\in\N}\Ext^i_A(A/I^n,M^\bullet)\isom
R^i\Gamma_{\!\!Z} M^{\bullet\sim}
\qquad\text{for every $i\in\N$}.
$$
\end{theorem}
\begin{proof} For $\cF^\bullet:=M^{\bullet\sim}$, trivial duality
(theorem \ref{th_trivial-dual}) identifies the source of
\eqref{eq_ext-and-gamma_Z} with $R\Hom^\bullet_A(A/I^n,M^\bullet)$;
then one takes cohomology in degree $i$ and forms the colimit over
$n$ to define the sought map. Next, by usual spectral sequence
arguments, we may reduce to the case where $M^\bullet$ is a single
$A$-module $M$ sitting in degree zero (see {\em e.g.} the proof of
proposition \ref{prop_der-cohereur}(i)).
By inspecting the definitions, one verifies easily that
the morphism thus defined is the composition of
the isomorphism of proposition \ref{prop_depth-Kosz}(iii)
and the map \eqref{eq_col-Exts} (see the proof of proposition
\ref{prop_depth-Kosz}(iii)); then it suffices to show
that the inverse system $(H_i\bK_\bullet(\bff^n)~|~n\in\N)$
is essentially zero when $i>0$ (lemma \ref{lem_ess-zero}).
By lemma \ref{lem_Hartsho}, this will in turn follow, provided
the following holds. For every finitely presented quotient $B$
of $A$, and every $b\in B$, there exists $p\in\N$ such that
$\Ann_B(b^q)=\Ann_B(b^p)$ for every $q\geq p$. This latter
assertion is a special case of corollary \ref{cor_prim-dec-exist}.
\end{proof}

\begin{corollary}\label{cor_calculate_Z}
In the situation of theorem {\em\ref{th_calculate_Z}}, we have
natural isomorphisms :
$$
\colim_{n\in\N}\Ext^i_A(I^n,M^\bullet)\isom
H^i(X\!\setminus\! Z,M^{\bullet\sim})
\qquad\text{for every $i\in\N$}.
$$
\end{corollary}
\begin{proof} Set $U:=X\!\setminus\! Z$. We may assume that $M^\bullet$
is a complex of injective $A$-modules, in which case the sought
map is obtained by taking colimits over the direct system
of composed morphisms :
$$
\Hom_A(I^n,M^\bullet)\xrightarrow{\beta_n}
\Gamma(U,M^\bullet)\to R\Gamma(U,M^{\bullet\sim})
$$
where $\beta_n$ is induced by the identification
$(I^n)^\sim_{|U}=\cO_{\!U}$ and the natural isomorphism :
$$
\Hom_{\cO_{\!U}}(\cO_{\!U},M^{\bullet\sim})\simeq
\Gamma(U,M^{\bullet\sim}).
$$
The complex $R\Gamma(U,M^{\bullet\sim})$ is computed by a
Cartan-Eilenberg injective resolution $M^{\bullet\sim}\isom\cM^\bullet$
of $\cO_{\!X}$-modules, and then the usual arguments allow to
reduce to the case where $M^\bullet$ consists of a single
injective $A$-module $M$ placed in degree zero.
Finally, from the short exact sequence $0\to I^n\to A\to A/I^n\to 0$
we deduce a commutative ladder with exact rows :
$$
\xymatrix{
0 \ar[r] & \Hom_A(A/I^n,M) \ar[r] \ar[d]_{\alpha_n} &
M \ar[r] \ddouble & \Hom_A(I^n,M) \ar[r] \ar[d]^{\beta_n} & 0 \ar[d] \\
0 \ar[r] & \Gamma_{\!\!Z}M^\sim \ar[r] & M \ar[r] &
\Gamma(U,M^\sim) \ar[r] & R^1\Gamma_{\!\!Z}M^\sim \ar[r] & 0
}$$
where $\alpha_n$ is induced from \eqref{eq_ext-and-gamma_Z}.
From theorem \ref{th_calculate_Z} it follows that
$\colim_{n\in\N}\alpha_n$ is an isomorphism, and
$R^i\Gamma_{\!\!Z}M^\sim$ vanishes for all $i>0$;
hence $\colim_{n\in\N}\beta_n$ is an isomorphism,
and the contention follows.
\end{proof}

\begin{proposition}\label{prop_duality-smooth}
Let $X\xrightarrow{f}Y\xrightarrow{g}S$ be two finitely
presented morphisms of schemes.
\begin{enumerate}
\item
If $g$ is smooth, $\cO_Y[0]$ is a dualizing complex on $Y$.
\item
If $g$ is smooth and $f$ is a closed immersion, then $X$ admits
a dualizing complex $\omega^\bullet_X$.
\item
If $X$ and $Y$ are affine and $\omega_Y^\bullet$ is a dualizing
complex on $Y$, then $f^!\omega_Y^\bullet$ is dualizing on $X$
(notation of \eqref{subsec_sharp-flat}).
\item
Every finitely presented quasi-separated $S$-scheme admits a
dualizing complex.
\end{enumerate}
\end{proposition}
\begin{proof} (i) is an immediate consequence of proposition
\ref{prop_hdim-finite}(i).

(ii): As in the proof of corollary \ref{cor_not-so-triv}, this
follows directly from (i) and lemma \ref{lem_transit-dual}(i).

(iii): Let us choose a factorization $f=p_Y\circ i$ where $i:X\to\A^n_Y$
is a finitely presented closed immersion, and $p_Y:\A^n_Y\to Y$
the smooth projection onto $Y$. In view of lemma \ref{lem_transit-dual}(i)
and \eqref{subsec_sharp-flat}, it suffices to prove the
assertion for the morphism $p_Y$. To this aim, we pick a closed finitely
presented immersion $h:Y\to\A^m_S$ and consider the fibre diagram :
$$
\xymatrix{ \A^n_Y \ar[r]^{h'} \ar[d]_{p_Y} & \A^{n+m}_S \ar[d]^{p_S} \\
Y \ar[r]^h & \A^m_S.
}$$
By (i) we know that the scheme $\A^m_S$ admits a dualizing
complex $\omega^\bullet$, and then lemma \ref{lem_transit-dual}(i)
says that $h^!\omega^\bullet$ is dualizing as well.
By proposition \ref{prop_unique-dual} we deduce that
$\omega^\bullet_Y\simeq\cL[\sigma]\otimes_{\cO_Y}h^!\omega^\bullet$
for some invertible $\cO_Y$-module $\cL$ and some continuous
function $\sigma:|Y|\to\Z$. Since $p_Y$ is smooth,
we can compute: $p^!_Y\omega^\bullet_Y\simeq
(p^!_Y\circ h^!\omega^\bullet)\otimes_{\cO_{\A^n_Y}}p_Y^*\cL[\sigma]$,
hence $p^!_Y\omega^\bullet_Y$ is dualizing if and only if
the same holds for $p^!_Y\circ h^!\omega^\bullet$. By proposition
\ref{prop_sharp-flat}(iv), the latter complex is isomorphic
to $h^{'\flat}\circ p_S^!\omega^\bullet$ and again using
lemma \ref{lem_transit-dual}(i) we reduce to checking that
$p_S^!\omega^\bullet$ is dualizing, which is clear from (i).

(iv): Let $f:X\to S$ be a finitely presented morphism, with $X$
quasi-separated. If $X$ is affine, (iii) says that $f^!\cO_{\!S}[0]$
is dualizing on $X$. In the general case, let $(U_i~|~i=1,\dots,n)$
be a finite covering of $X$ consisting of affine open subsets;
for each $i,j,k=1,\dots,n$, denote by $f_i:U_i\to S$ the restriction
of $f$, set $U_{ij}:=U_i\cap U_j$, $U_{ijk}:=U_{ij}\cap U_k$, and
let $g_{ij}:U_{ij}\to U_i$ be the inclusion map. We know that
$f^!_i\cO_{\!S}[0]$ is dualizing on $U_i$, for every $i=1,\dots,n$;
moreover, for every $i,j=1,\dots,n$ there exists a natural isomorphism
$$
\psi_{ij}:g^*_{ij}f^!_i\cO_{\!S}[0]\isom g^*_{ji}f^!_j\cO_{\!S}[0]
$$
fulfilling the cocycle condition
$$
\psi_{jk|U_{ijk}}\circ\psi_{ij|U_{ijk}}=\psi_{ik|U_{ijk}}
\qquad
\text{for every $i,j,k=1,\dots,n$}
$$
(lemma \ref{lem_compat-for_!}). In other words,
$((U_i,f^!_i\cO_{\!S}[0]);\psi_{ij}~|~i,j=1,\dots,n)$ is
a descent datum for the fibration \eqref{eq_dualiz-fibrat}.
Then the assertion follows from proposition \ref{prop_glue-dualizers}.
\end{proof}

\begin{example}\label{ex_duality-smooth}
(i)\ \
For given $b\in\fm_K$, let $i_b:S_{\!/b}\to S$ be the closed
immersion. If $b\neq 0$, a simple computation yields a natural
isomorphism in $\sD(\cO_{\!S_{\!/b}}\Mod)$ :
$$
i^!_b\cO_{\!S}\isom\cO_{\!S_{\!/b}}[-1].
$$
By proposition \ref{prop_duality-smooth}(i,iii), we deduce
that $\cO_{\!S_{\!/b}}[0]$ is a dualizing complex on $S_{\!/b}$,
for every $b\in\fm_K$.

(ii)\ \
Next, let $f:X\to S_{\!/b}$ be an affine finitely presented
Cohen-Macaulay morphism, of constant fibre dimension $n$.
Then $\omega^\bullet_X:=f^!\cO_{S_{\!/b}}[0]$ is a dualizing
complex on $X$, by (i) and proposition \ref{prop_duality-smooth}(iii).
Moreover, $\omega^\bullet_X$ is concentrated in degree $-n$,
and $H^{-n}\omega^\bullet_X$ is a finitely presented
$f$-Cohen-Macaulay $\cO_{\!X}$-module.
Indeed, let $i:X\to Y:=\A^m_{S_{\!/b}}$ be a closed immersion
of $S_{\!/b}$-schemes, and denote by $g:Y\to S_{\!/b}$ the
projection. Fix any $x\in X$, let $y:=i(x)$ and set
$$
d_x:=\dim\cO_{f^{-1}(fx),x}
\qquad
d_y:=\dim\cO_{g^{-1}(gy),y}.
$$
We have a natural isomorphism in $\sD(\cO_{\!X}\Mod)$
$$
\omega^\bullet_X\isom i^!\cO_Y[m]\isom
i^*R\cHom^\bullet_{\cO_Y}(i_*\cO_{\!X},\cO_Y[m])
$$
(lemma \ref{lem_compat-for_!}(i)), and by assumption
$i_*\cO_{\!X}$ is a $g$-Cohen-Macaulay $\cO_Y$-module;
by theorem \ref{th_CM-duality}(ii), it follows that
$\omega^\bullet_{X,x}$ is concentrated in degree
$d_y-d_x-m$, and its cohomology in that degree is
$f$-Cohen-Macaulay, as asserted. Lastly, since the
fibres of $f$ and $g$ are equidimensional
(\cite[Ch.IV, Prop.5.2.1]{EGAIV-2} and lemma
\ref{lem_CM}(ii)), it is easily seen that $d_y-d_x=m-n$
(details left to the reader), whence the claim. 
\end{example}

\sset\subsubsection{}\label{subsec_def-cod-map}
For any finitely presented morphism $f:X\to S$ we consider the map :
$$
d:|X|\to\Z\qquad
x\mapsto\tr.\deg(\kappa(x)/\kappa(f(x)))+\dim\overline{\{f(x)\}}.
$$
\begin{lemma}\label{lem_cod-function}
With the notation of \eqref{subsec_def-cod-map}, let $x,y\in X$,
and suppose that $x$ is a specialization of $y$. We have :
\begin{enumerate}
\item
$x$ is an immediate specialization of $y$ if and only if
$d(y)=d(x)+1$.
\item
If $X$ is irreducible, $d(x)-d(y)=\dim X(y)-\dim X(x)$.
\item
$X$ is catenarian and of finite Krull dimension.
\item
If $X$ is irreducible and $f$ is flat, then $\cO_{\!f^{-1}(fx),x}$
is equidimensional.
\item
If $f$ is flat, then
$\dim\cO_{\!X,x}=\dim\cO_{\!f^{-1}(fx),x}+\dim\cO_{\!S,f(x)}$.
\item
If $f$ is Cohen-Macaulay at the point $x$, then $\cO_{\!X,x}$
is equidimensional.
\end{enumerate}
\end{lemma}
\begin{proof}(i): In case $f(x)=f(y)$, the assertion follows from
\cite[Ch.IV, Prop.5.2.1]{EGAIV-2}. Hence, suppose that $f(x)=s$,
$f(y)=\eta$; let $Z$ be the topological closure of $\{y\}$ in
$X$, and endow $Z$ with its reduced subscheme structure; notice
that $Z$ is an $S$-scheme of finite type. By assumption, $X$ is
quasi-compact and quasi-separated, hence $\{y\}$ is a
pro-constructible subset of $X$, and therefore $Z$ is the set
of all specializations of $y$ in $X$
(\cite[Ch.IV, Th.1.10.1]{EGAIV}). Especially, $y$ is the unique
maximal point of $Z_\eta:=Z\cap f^{-1}(\eta)$; also, if $x$ is
an immediate specialization of $y$, then $x$ is a maximal point
of $Z_s:=Z\cap f^{-1}(s)$, and by \cite[Ch.IV, Lemme 14.3.10]{EGAIV-3},
the latter implies that $d(y)=d(x)+1$.

Conversely, suppose that $d(y)=d(x)+1$; from
\cite[Ch.IV, Lemme 14.3.10]{EGAIV-3} and
\cite[Ch.IV, Prop.5.2.1]{EGAIV-2} we deduce that $x$ is a maximal
point of $Z_s$; then $x$ is an immediate specialization of $y$ in
$X$, since otherwise $x$ would be a specialization in $X$
of a proper specialization $y'$ of $y$ in $Z_\eta$, and in
this case \cite[Ch.IV, Lemme 14.3.10]{EGAIV-3} would say that
$d(y')=d(x)+1$, {\em i.e.} $d(y)=d(y')$, contradicting
\cite[Ch.IV, Prop.5.2.1]{EGAIV-2}.

(iii) It is easily seen that $X$ has finite Krull dimension.
Now, consider any sequence $y_0,\dots,y_n$ of points
of $X$, with $y_0:=y$, $y_n:=x$, and such that $y_{i+1}$ is an
immediate specialization of $y_i$, for $i=0,\dots,n-1$; from
(i) we deduce that $n=d(y)-d(x)$, especially $n$ is independent
of the chosen chain of specializations, so $X$ is catenarian.

(ii): We reduce easily to the case where $y$ is the maximal
point of $X$, in which case we may argue as in the proof of
(iii) (details left to the reader).

(iv): This is trivial, if $f(x)=\eta$; hence let us assume
that $f(x)=s$, and let $y_1,y_2$ be two maximal generizations
of $x$ in $f^{-1}(s)$. We may find an affine open subset
$U\subset X$ such that $y_1\in U$ and $y_2\notin U$; arguing
as in the proof of (i), we see that the maximal point $z$ of
$U$ (which is also the maximal point of $X$) is an immediate
generization of $y_1$ in $U$, hence $d(z)=d(y_1)+1$. Likewise,
$d(z)=d(y_2)+1$, hence $d(y_1)=d(y_2)$. In view of (i), the
assertion follows easily.

(v): If $f(x)=\eta$, the identity is \cite[Ch.IV, Cor.6.1.2]{EGAIV-2}.
Hence, suppose that $f(x)=s$. In this case, the proof of (iv)
shows that the identity holds, when $X$ is irreducible.
In the general case, notice that -- by the flatness assumption --
every irreducible component of $X$ intersect $f^{-1}(\eta)$,
and it is therefore itself a flat $S$-scheme (with its reduced
subscheme structure). Since $f^{-1}(\eta)$ is a noetherian
scheme, it also follows that the set of irreducible components
of $X$ is finite; thus, let $X_1,\dots,X_n$ be the reduced
irreducible components of $X$ containing $x$, and $f_i:X_i\to S$
($i=1,\dots,n$) the corresponding restrictions of $f$. Then
$$
\dim\cO_{\!X,x}=\max_{1\leq i\leq n}\dim\cO_{\!X_i,x}
\qquad\text{and}\qquad
\dim\cO_{\!f^{-1}(s),x}=\max_{1\leq i\leq n}\dim\cO_{\!f_i^{-1}(s),x}
$$
whence the contention.

(vi): By lemma \ref{lem_CM}(ii), the morphism $f$ is equidimensional
at the point $x$. Especially, $f(u)=\eta$ for every maximal
point $u$ of $\Spec\,\cO_{\!X,x}$, and there exist an open
neighborhood $U$ of $x$ in $X$ and an integer $e\in\N$,
such that $f^{-1}f(y)$ is equidimensional of dimension $e$,
for every $y\in U$. Since $f$ is finitely presented, it follows
that
$$
\tr.\deg(\kappa(y)/\kappa(f(y)))=e
$$
for every $y\in U$ such that $y$ is maximal in $f^{-1}f(y)$
(\cite[Th.5.6]{Mat}). Hence, $d(u)=e+1$ for every maximal
point $u$ of $\Spec\,\cO_{\!X,x}$, and especially, for every
such $u$, the difference $d(x)-d(u)$ is independent of $u$;
in light of (i), the assertion follows.
\end{proof}

\begin{lemma}\label{lem_coh-inject}
{\em (i)}\ \
The $K^+$-module $K/K^+$ is coh-injective.
\begin{enumerate}
\addenu
\item
Let $A$ be a finitely presented $K^+$-algebra, $I\subset A$
a finitely generated ideal, $J$ a coh-injective $K^+$-module.
Then the $A$-module
$$
J_A:=\colim_{n\in\N}\Hom_{K^+}(A/I^n,J)
$$
is coh-injective.
\end{enumerate}
\end{lemma}
\begin{proof} (i): It suffices to show that $\Ext^1_{K^+}(M,K/K^+)=0$
whenever $M$ is a finitely presented $K^+$-module.
This is clear when $M=K^+$, and then \cite[Ch.6, Lemma 6.1.14]{Ga-Ra}
reduces to the case where $M=K^+/aK^+$ for some
$a\in\fm_K\!\setminus\!\{0\}$; in that case we can compute
using the free resolution $K^+\xrightarrow{a}K^+\to M$,
and the claim follows easily.

(ii): According to theorem \ref{th_omega-coh}(ii.a), the
$A$-module $J'_A:=\Hom_{K^+}(A,J)$ is coh-injective. However,
$J_A=\bigcup_{n\in\N}\Ann_{J'_A}(I^n)$, so the assertion follows
from theorem \ref{th_omega-coh}(ii.b).
\end{proof}

\begin{theorem}\label{th_stalk-of-omega}
Let $f:X\to S$ be a finitely presented affine morphism.
Then for every point $x\in X$, the $\cO_{\!X,x}$-module
$$
J(x):=R^{1-d(x)}\Gamma_{\!\{x\}}(f^!\cO_{\!S}[0])_{|X(x)}
$$
is coh-injective, and we have a natural isomorphism
in $\sD(\cO_{X,x}\Mod)$ :
$$
R\Gamma_{\!\{x\}}(f^!\cO_{\!S}[0])_{|X(x)}\simeq J(x)[d(x)-1].
$$
\end{theorem}
\begin{proof} Fix $x\in X$ and set $d:=d(x)$.

\begin{claim} The theorem holds in case $f(x)=\eta$.
\end{claim}
\begin{pfclaim} Indeed, let
$f_\eta:f^{-1}(\eta)\to S(\eta):=\Spec\,K$ be the restriction
of $f$; according to lemma \ref{lem_compat-for_!}(i),
$f^!_\eta\cO_{\!S(\eta)}[0]$ is the restriction of $f^!\cO_{\!S}[0]$,
so the claim follows immediately from example \ref{ex_compute-c_X}.
\end{pfclaim}

Hence, suppose that $f(x)=s$, the closed point of $S$, and say
that $X=\Spec\,A$.

\begin{claim}\label{cl_reduce-to-cl}
We may assume that $x$ is closed in $X$, hence that $\kappa(x)$
is finite over $\kappa(s)$.
\end{claim}
\begin{pfclaim} Arguing as in the proof of lemma
\ref{lem_integral-fp-ext}, we can find a factorization
of the morphism $f$ as a composition
$X\xrightarrow{g}Y:=\A^d_S\xrightarrow{h}S$, such that
$\xi:=g(x)$ is the generic point of $h^{-1}(s)\subset Y$,
the morphism $g$ is finitely presented, and the stalk
$\cO_{Y,\xi}$ is a valuation ring. Moreover, let
$g_y:=g\times_Y\one_{Y(y)}:X(y):=X\times_YY(y)\to Y(y)$;
we have a natural isomorphism
$$
(f^!\cO_{\!S}[0])_{|X(y)}\simeq g_y^!\cO_{Y(y)}[d].
$$
Then we may replace $f$ by $g_y$, and $K^+$ by $\cO_{Y,\xi}$,
whence the claim.
\end{pfclaim}

Hence, suppose now that $x$ is closed in $X$, let $\fp\subset A$
be the maximal ideal corresponding to $x$,
and $\bar\fp$ the image of $\fp$ in $\bar A:=A\otimes_{K^+}\kappa$;
we choose a finite system of elements $b_1,\dots,b_r\in\cO_{\!X}(X)$
whose images in $\bar A$ generate $\bar\fp$. Pick any non-zero
$a\in\fm_K$, and let $I\subset A$ be the ideal generated by the
system $(a,b_1,\dots,b_r)$, and $\cI\subset\cO_{\!X}$ the corresponding
coherent ideal; clearly $V(I)=\{x\}$, hence theorem
\ref{th_calculate_Z} yields a natural isomorphism :
\set\begin{equation}\label{eq_colim_Z}
\colim_{n\in\N}R^i\Hom^\bullet_{\cO_{\!X}}(\cO_{\!X}/\cI^n,f^!\cO_{\!S}[0])
\isom R^i\Gamma_{\!\{x\}} f^!\cO_{\!S}[0]
\qquad\text{for every $i\in\Z$}
\end{equation}
where the transition maps in the colimit are induced by
the natural maps $\cO_{\!X}/\cI^n\to\cO_{\!X}/\cI^m$, for
every $n\geq m$.
However, from lemmata \ref{lem_rightist}(ii), \ref{lem_compat-for_!}(i)
and corollary \ref{cor_Ext-loc=glob}(i) we deduce as well natural
isomorphisms :
$$
R\Hom^\bullet_{\cO_{\!X}}(\cO_{\!X}/\cI^n,f^!\cO_{\!S}[0])\isom
R\Hom^\bullet_{\cO_{\!S}}(f_*\cO_{\!X}/\cI^n,\cO_{\!S})\isom
R\Hom^\bullet_{K^+}(A/I^n,K^+).
$$
We wish to compute these $\Ext$ groups by means of lemma
\ref{lem_omega-coh}(iii); to this aim, let us remark first
that $A/I^n$ is an integral $K^+$-algebra for every $n\in\N$,
hence it is a finitely presented torsion $K^+$-module,
according to proposition \ref{prop_integral-fp-ext}(i). We may
then use the coh-injective resolution $K^+[0]\to(0\to K\to K/K^+\to 0)$
(lemmata \ref{lem_coh-inject}(i), \ref{lem_omega-coh}(iii)) to compute :
$$
\Ext^i_{K^+}(A/I^n,K^+)=\left\{
\begin{array}{ll} J_n:=\Hom_{K^+}(A/I^n,K/K^+) & \text{if $i=1$} \\
                  0 & \text{otherwise.}
\end{array}\right.
$$
Now the theorem follows from lemma \ref{lem_coh-inject}(i)
and theorem \ref{th_omega-coh}(ii.a),(iii).

We could also appeal directly to lemma \ref{lem_coh-inject}(ii),
provided we already knew that the foregoing natural identifications
transform the direct system whose colimit appears in \eqref{eq_colim_Z},
into the direct system $(J_n~|~n\in\N)$ whose transition maps are
induced by the natural maps $A/I^n\to A/I^m$, for every $n\geq m$.
For the sake of completeness, we check this latter assertion.

For every $n\geq m$, let $j_n:X_n:=\Spec\,A/I^n\to X$ and
$j_{mn}:X_m\to X_n$ be the natural closed immersions.
We have a diagram of functors :
$$
\xymatrix{ j_{m*}\circ j^!_{mn}\circ(f\circ j_n)^!
\ar[rr]^-{\zeta_1} \ar[rrd]^{\eps_1} & &
j_{n*}\circ(j_{mn*}\circ j_{mn}^!)\circ j_n^!\circ f^!
\ar[rrd]^{\eps_2} & &
j_{m*}\circ j_m^!\circ f^! \ar[ll]_-{\xi_2} \ar[d]^{\alpha_{mn}} \\
j_{m*}\circ(f\circ j_m)^! \ar[u]^{\xi_1} \ar[rr]^-{\beta_{mn}} & &
j_{n*}\circ(f\circ j_n)^! \ar[rr]^-{\zeta_2} & & j_{n*}\circ j_n^!\circ f^!
}$$
where :
\begin{itemize}
\item
$\zeta_1$ and $\zeta_2$ are induced by the natural isomorphism
$\psi_{f,j_n}:(f\circ j_n)^!\isom j_n^!\circ f^!$ of lemma
\ref{lem_compat-for_!}(i).
\item
$\xi_1$ (resp. $\xi_2$) is induced by the isomorphism
$\psi_{f\circ j_n,j_{mn}}$ (resp. $\psi_{j_n,j_{mn}}$).
\item
$\eps_1$ and $\eps_2$ are induced by the counit of adjunction
$j_{mn*}\circ j_{mn}^!\to\one_{\sD^+(\cO_{\!X_n}\Mod)}$.
\item
$\beta_{mn}$ (resp. $\alpha_{mn}$) is induced by the natural map
$(f\circ j_n)_*\cO_{\!X_n}\to(f\circ j_m)_*\cO_{\!X_m}$ (resp.
by the map $j_{n*}\cO_{\!X_n}\to j_{m*}\cO_{\!X_m}$).
\end{itemize}
It follows from lemma \ref{lem_wants-to-commute} that the two
triangular subdiagrams commute, and it is also clear that the
same holds for the inner quadrangular subdiagram.
Moreover, lemma \ref{lem_compat-for_!}(ii) yields a commutative
diagram
\set\begin{equation}\label{eq_bim-bum}
{\diagram (f\circ j_m)^!
\ar[d]_{\psi_{f\circ j_n,j_{mn}}} \ar[rr]^{\psi_{f,j_m}} & &
j_m^!\circ f^! \ar[d]^{\psi_{j_n,j_{mn}}\circ f^!} \\
j^!_{mn}\circ(f\circ j_n)^! \ar[rr]^{j^!_{mn}(\psi_{f,j_n})} & &
j_{mn}^!\circ j_n^!\circ f^!
\enddiagram}\end{equation}
such that $j_{m*}\eqref{eq_bim-bum}$ is the diagram :
$$
\xymatrix{ j_{m*}\circ(f\circ j_m)^! \ar[d]_{\xi_1}
\ar[rr]^{j_{m*}(\psi_{f,j_m})} & &
j_{m*}\circ j_m^!\circ f^! \ar[d]^{\xi_2} \\
j_{m*}\circ j^!_{mn}\circ(f\circ j_n)^! \ar[rr]^{\zeta_1} & &
j_{m*}\circ j_{mn}^!\circ j_n^!\circ f^!.
}$$
We then arrive at the following commutative diagram :
$$
\xymatrix{
R\Hom^\bullet_{\cO_{\!S}}(\cO_{\!X}/\cI^m,\cO_{\!S})
\ar[rr]^-{R\Gamma(\psi_{f,j_m})} \ar[d]_{R\Gamma(\beta_{mn})} & &
R\Hom^\bullet_{\cO_{\!X}}(\cO_{\!X}/\cI^m,f^!\cO_{\!S}[0])
\ar[d]^{R\Gamma(\alpha_{mn})} \\
R\Hom^\bullet_{\cO_{\!S}}(\cO_{\!X}/\cI^n,\cO_{\!S})
\ar[rr]^-{R\Gamma(\psi_{f,j_n})} & &
R\Hom^\bullet_{\cO_{\!X}}(\cO_{\!X}/\cI^n,f^!\cO_{\!S}[0]).
}$$
Now, the maps $R\Gamma(\alpha_{mn})$ are the transition morphisms
of the inductive system whose colimit appears in \eqref{eq_colim_Z},
hence we may replace the latter by the inductive system formed by
the maps $R\Gamma(\beta_{mn})$. Combining with corollary
\ref{cor_Ext-loc=glob}(i), we finally deduce natural $A$-linear
isomorphisms
$$
\colim_{n\in\N}\Ext^i_{K^+}(A/I^n,K^+)\isom
R^i\Gamma_{\!\{x\}} f^!\cO_{\!S}[0]
$$
where the transition maps in the colimit are induced by
the natural maps $A/I^n\to A/I^m$, for every $n\geq m$.
Our assertion is an immediate consequence.
\end{proof}

\sset\subsubsection{}\label{subsec_loc-duality}
Let $A$ be a local ring, set $X:=\Spec\,A$ let $x\in X$
be the closed point, and suppose that
\begin{enumerate}
\alphaenu
\item
either, $A$ is noetherian and $X$ admits a dualizing complex
\item
or else, $A$ is essentially of finite presentation over $K^+$.
\end{enumerate}
Notice that $X$ admits a dualizing complex in case (b) as well.
Indeed, in that case one can find a finitely presented affine
$S$-scheme $Y$ and a point $y\in Y$ such that
$X\simeq\Spec\,\cO_{Y,y}$; by proposition
\ref{prop_duality-smooth}(iv), $Y$ admits a dualizing complex
$\omega^\bullet_Y$, and since every coherent $\cO_{\!X}$-module
extends to a coherent $\cO_Y$-module, one verifies easily that
the restriction of $\omega^\bullet_Y$ is dualizing for $X$ (cp.
the proof of proposition \ref{prop_fishy}).
Hence, let $\omega^\bullet$ be a dualizing complex for $X$.
In case (a) (resp. in case (b)), it follows easily from
corollary \ref{cor_dual-catenary} (resp. from theorem
\ref{th_stalk-of-omega} and propositions
\ref{prop_duality-smooth}(iii) and \ref{prop_unique-dual},
that there exists a unique $c\in\Z$ such that
$$
J(x):=R^c\Gamma_{\!\{x\}}\omega^\bullet\neq 0.
$$
Moreover, $J(x)$ is a coh-injective $A$-module, hence by lemma
\ref{lem_omega-coh}(iii), we obtain a well defined functor
\set\begin{equation}\label{eq_dual-on-A-mod}
\sD^b(A\Mod_\coh)\to\sD^b(A\Mod)^o
\quad :\quad C^\bullet\mapsto\Hom^\bullet_A(C^\bullet,J(x)).
\end{equation}
Furthermore, let $\sD^b_{\{x\}}(A\Mod_\coh)$ be the full
subcategory of $\sD^b(A\Mod_\coh)$ consisting of all complexes
$C^\bullet$ such that $\Supp\,H^\bullet C^\bullet\subset\{x\}$.
We have the following :

\begin{corollary}[Local duality]\label{cor_loc-duality}
In the situation of \eqref{subsec_loc-duality}, the following
holds :
\begin{enumerate}
\item
The functor \eqref{eq_dual-on-A-mod} restricts to an
equivalence of categories :
$$
D:\sD^b_{\{x\}}(A\Mod_\coh)\to\sD^b_{\{x\}}(A\Mod_\coh)^o
$$
and the natural transformation $C^\bullet\to D\circ D(C^\bullet)$
is an isomorphism of functors.
\item
For every $i\in\Z$ there exists a natural isomorphism of functors :
$$
R^i\Gamma_{\!\{x\}}\circ\cD\isom D\circ R^{c-i}\Gamma
\quad : \quad \sD^b(\cO_{\!X}\Mod)_\coh\to A\Mod^o
$$
where $\cD$ is the duality functor corresponding to $\omega^\bullet$.
\item
Let $c$ and $\cD$ be as in {\em(ii)}, $I\subset A$ a finitely
generated ideal such that $V(I)=\{x\}$, and denote by $A^\wedge$
the $I$-adic completion of $A$. Then for every $i\in\Z$ there
exists a natural isomorphism of functors :
$$
D\circ R^i\Gamma_{\!\{x\}}\circ\cD\isom
A^\wedge\otimes_A R^{c-i}\Gamma
\quad : \quad \sD^b(\cO_{\!X}\Mod)_\coh\to A^\wedge\Mod_\coh^o.
$$
\end{enumerate}
\end{corollary}
\begin{proof} Let $C_\bullet$ be any object of
$\sD^b_{\{x\}}(A\Mod_\coh)$; we denote by $i:\{x\}\to X$
the natural closed immersion. Obviously
$C^\sim_\bullet=i_*i^{-1}C^\sim_\bullet$ (notation of
\eqref{sec_various-O-mod}), therefore :
$$
\begin{aligned}
\Hom^\bullet_A(C^\bullet,J(x)) & \isom R\Hom^\bullet_A(C^\bullet,J(x))
& \quad & \text{by lemma \ref{lem_omega-coh}(iii)} \\
& \isom R\Hom^\bullet_{\cO_{\!X}}(C^\sim_\bullet,
R\underline\Gamma_{\{x\}}\omega^\bullet[c])
& \quad & \text{by theorem \ref{th_trivial-dual}} \\
& \isom R\Hom^\bullet_{\cO_{\!X}}(i_*i^{-1}C^\sim_\bullet,\omega^\bullet[c])
& \quad & \text{by lemma \ref{lem_adj-Gamma_Z}(i.b)} \\
& \isom R\Hom^\bullet_{\cO_{\!X}}(C^\sim_\bullet,\omega^\bullet[c])
\end{aligned}
$$
which easily implies (i). Next, we compute, for any object
$\cF^\bullet$ of $\sD^b(\cO_{\!X}\Mod)_\coh$ :
$$
\begin{aligned}
R^i\Gamma_{\!\{x\}}\cD(\cF^\bullet) & \isom
R^i\Hom^\bullet_{\cO_{\!X}}(\cF^\bullet,R\Gamma_{\!\{x\}}\omega^\bullet) &
\quad & \text{by lemma \ref{lem_adj-Gamma_Z}(iii)} \\
& \isom R^i\Hom^\bullet_{\cO_{\!X}}(\cF^\bullet,J(x)[-c]) \\
& \isom R^i\Hom^\bullet_A(R\Gamma\cF^\bullet,J(x)[-c]) &
\quad & \text{by theorem \ref{th_trivial-dual}} \\
& \isom \Hom_A(R^{c-i}\Gamma\cF^\bullet,J(x)) &
\quad & \text{by lemma \ref{lem_omega-coh}(iii)}
\end{aligned}
$$
whence (ii). Finally, let $\cF^\bullet$ be any
object of $\sD^b(\cO_{\!X}\Mod)_\coh$; if $A$ is as in case
(b) of \eqref{subsec_loc-duality}, we compute :
$$
\begin{aligned}
A^\wedge\otimes_AR^{c-i}\Gamma\cF^\bullet & \isom
\lim_{n\in\N}\,(A/I^n)\otimes_AR^{c-i}\Gamma\cF^\bullet
& \quad & \text{by corollary \ref{cor_cption-exact}(ii)} \\
& \isom\lim_{n\in\N}\,
D\circ D((A/I^n)\otimes_AR^{c-i}\Gamma\cF^\bullet)
& \quad & \text{by (i)}\\
& \isom\lim_{n\in\N}\,
D(\Hom_A(A/I^n,R^i\Gamma_{\!\{x\}}\cD(\cF^\bullet)))
& \quad & \text{by (ii)} \\
& \isom D(\colim_{n\in\N}\Hom_A(A/I^n,R^i\Gamma_{\!\{x\}}\cD(\cF^\bullet))) \\
& \isom D(R^i\Gamma_{\!\{x\}}\cD(\cF^\bullet))
\end{aligned}
$$
so (iii) holds as well. If $A$ is as in case (a) of
\eqref{subsec_loc-duality}, the same proof works : instead
of corollary \ref{cor_cption-exact}(ii), one has just to
appeal to \cite[Th.8.7 and Th.8.8]{Mat}.
\end{proof}

\begin{definition} Let $A$ and $J(x)$ be as in
\eqref{subsec_loc-duality}, and $E$ an $A$-module. We say that $E$
is {\em finitely copresented} if there exist a finitely presented
$A$-module $M$ and an $A$-linear isomorphism $E\isom\Hom_A(M,J(x))$.
\end{definition}

\begin{lemma}\label{lem_copresentation}
Let $A$ be as in \eqref{subsec_loc-duality}, and $I\subset A$ a
finitely generated ideal, whose radical is the maximal ideal of
$A$. For every $A$-module $E$, the following conditions are
equivalent :
\begin{enumerate}
\alphaenu
\item
$E$ is finitely copresented and finitely generated.
\item
$E$ is finitely copresented and $I^kE=0$ for some $k\in\N$.
\item
$E$ is finitely presented and $I^kE=0$ for some $k\in\N$.
\end{enumerate}
\end{lemma}
\begin{proof} To ease notation, for every ideal $L\subset A$
and every $A$-module $M$ set 
$$
D(M):=\Hom_A(M,J(x))
\quad\text{and}\quad
M[L]:=\Ann_M(L).
$$

\begin{claim}\label{cl_copresentation}
Let $E$ be any finitely copresented $A$-module. Then
$E[I^k]$ is a finitely presented $A$-module for every
$k\in\N$, and $E=\bigcup_{k\in\N}E[I^k]$.
\end{claim}
\begin{pfclaim} We may assume that $E=D(M)$ for some
finitely presented $A$-module $M$; then clearly every
element of $E$ is annihilated by some power of $I$, since
the same holds for $J(x)$. Moreover, we have a natural
identification
$$
E[I^n]\isom D(M/I^nM)
\qquad
\text{for every $n\in\N$}.
$$
But $M/I^kM$ is finitely presented and supported at $\{x\}$,
so $D(M/I^kM)$ is finitely presented as well, by corollary
\ref{cor_loc-duality}(i).
\end{pfclaim}

It follows immediately from claim \ref{cl_copresentation}
that (a)$\Rightarrow$(b)$\Rightarrow$(c).

(c)$\Rightarrow$(a): If $E$ is finitely presented and
$I^kE=0$, then the same holds for $D(E)$, by corollary
\ref{cor_loc-duality}(i), and moreover $D\circ D(E)=E$,
by the same token, whence (a).
\end{proof}

\begin{proposition}\label{prop_copresentation}
Let $A$ be as in \eqref{subsec_loc-duality}, and $f:E\to N$
a homomorphism of $A$-modules, with $E$ finitely copresented
and $N$ finitely presented. Then we have :
\begin{enumerate}
\item
$\Ker\,f$ is a finitely copresented $A$-module.
\item
$\Img\,f$ and $\Coker\,f$ are finitely presented $A$-modules.
\end{enumerate} 
\end{proposition}
\begin{proof} Clearly, it suffices to show the assertions
concerning $\Ker\,f$ and $\Img\,f$. To this aim, pick any
finitely generated ideal $I\subset A$ such that
$\Supp\,A/I=\{x\}$ (where $x$ denotes the closed point of
$\Spec\,A$); by corollary \ref{cor_prim-dec-exist}(ii) there
exists $n\in\N$ such that $N[I^n]=N[I^m]$ for every integer
$m\geq n$. In view of claim \ref{cl_copresentation}, we see
that there exists $n\in\N$ such that $\Img\,f=(\Img\,f)[I^n]$.
Taking into account corollary \ref{cor_prim-dec-exist}(i), we
may then replace $N$ by $N[I^n]$, and assume from start that
$I^nN=0$, so $N$ is also finitely copresented, by lemma
\ref{lem_copresentation}.

Now, more generally, let $E$ and $F$ be two finitely copresented
$A$-modules, and say that $E=D(M)$, $F=D(P)$, where $D$ is defined
as in the proof of lemma \ref{lem_copresentation}; in view of claim
\ref{cl_copresentation} and corollary \ref{cor_loc-duality}(i)
we get natural identifications :
$$
\begin{aligned}
\Hom_A(E,F)\isom\, & \lim_{n\in\N}\Hom_A(E[I^n],F[I^n]) \\
\isom\, & \lim_{n\in\N}\Hom_A(P/I^nP,M/I^nM) \\
\isom\, & \lim_{n\in\N}\Hom_A(P,M/I^nM) \\
\isom\, & \lim_{n\in\N}\Hom_A(P,M^\wedge)
\end{aligned}
$$
where $M^\wedge$ denotes the $I$-adic completion of $M$.
Recall that the natural map $M\to M^\wedge$ is injective
(corollary \ref{cor_cption-exact}(ii,iii)), so the same
holds for the induced map
\set\begin{equation}\label{eq_induced-homs}
\Hom_A(P,M)\to\Hom_A(E,F)
\qquad
f\mapsto D(f).
\end{equation}

\begin{claim}\label{cl_both-copresent}
With the foregoing notation, we have :
\begin{enumerate}
\item
Let $g:E\to F$ be any morphism. If $g$ lies in the image of
\eqref{eq_induced-homs}, then $\Ker\,g$, $\Coker\,g$ and
$\Img\,g$ are all finitely copresented $A$-modules.
\item
If either $I^kM=0$ or $I^kP=0$ for some $k\in\N$, then
\eqref{eq_induced-homs} is an isomorphism.
\end{enumerate}
\end{claim}
\begin{pfclaim}(i): Indeed, say that $g$ is the image of
$h:P\to M$ under \eqref{eq_induced-homs}; since $J(x)$
is coh-injective, we get natural identifications
$$
\Ker\,g\isom D(\Coker\,h)
\qquad
\Coker\,g\isom D(\Ker\,h)
\qquad
\Img\,g\isom D(\Img\,h)
$$
whence the contention.

(ii): If $I^kM=0$ for some $k\in\N$, then the natural
map $M\to M^\wedge$ is an isomorphism, and the assertion
follows. If $I^kP=0$, then any morphism $P\to M^\wedge$
factors through $M^\wedge[I^k]$ (notation of the proof
of lemma \ref{lem_copresentation}); however, corollaries
\ref{cor_prim-dec-exist}(i) and \ref{cor_cption-exact}(i)
yield natural identifications
$$
M[I^n]\isom(M[I^n])^\wedge\isom M^\wedge[I^n]
$$
whence the contention (details left to the reader).
\end{pfclaim}

Now, letting $F:=N$ in claim \ref{cl_both-copresent}, we
deduce that both $\Ker\,f$ and $\Img\,f$ are finitely
copresented. Lastly, since $I^n\cdot\Img\,f=0$, lemma
\ref{lem_copresentation} says that $\Img\,f$ is
finitely presented, as sought.
\end{proof}

\begin{theorem}\label{th_finit}
In the situation of \eqref{subsec_loc-duality}, let
$I\subset A$ be any finitely generated ideal, $C^\bullet$
any object of $\sD^b(\cO_{\!X}\Mod_\coh)$, and set
$U:=X\!\setminus\!\{x\}$. For every $i\in\Z$, the following
conditions are equivalent :
\begin{enumerate}
\alphaenu
\item
$I\cdot H^i(U,C^\bullet)$ is a finitely presented $A$-module.
\item
$I\cdot H^i(U,C^\bullet)$ is a finitely generated $A$-module.
\item
$I\cdot R^{i+1}\Gamma_{\!\{x\}}C^\bullet$ is a finitely
presented $A$-module.
\item
 $I\cdot R^{i+1}\Gamma_{\!\{x\}}C^\bullet$ is a finitely
generated $A$-module.
\item
The radical of $\Ann_A(I\cdot R^{i+1}\Gamma_{\!\{x\}}C^\bullet)$
contains the maximal ideal of $A$.
\item
$I\cdot R^i\Gamma_{\!\{y\}}C^\bullet_{|X(y)}=0$ for every closed
point $y\in U$.
\end{enumerate}
\end{theorem}
\begin{proof} Let $\omega^\bullet$ be a dualizing complex for $X$
and $y\in X$ any point; then $\omega^\bullet_{|X(y)}$ is a dualizing
complex for $X(y)$, and according to \eqref{subsec_loc-duality},
there is a unique integer $c(y)$ such that
$J(y):=R^{c(y)}\Gamma_{\!\{y\}}\omega^\bullet_{|X(y)}\neq 0$.
If $A$ is as in case (b) (resp. as in case (a)) of
\eqref{subsec_loc-duality}, we may invoke theorem
\ref{th_stalk-of-omega} and lemma \ref{lem_cod-function}(i)
(resp. lemma \ref{lem_dual-catenary}) to see that
\set\begin{equation}\label{eq_tough-noam}
c(y)+1=c:=c(x)
\quad\Leftrightarrow\quad
\text{$y$ is a closed point of $U$}.
\end{equation}
By corollary \ref{cor_loc-duality}(i,ii), the rule
$$
M\mapsto D_y(M):=\Hom_A(M,J(y))
\qquad
\text{for every $\cO_{\!X,y}$-module $M$}
$$
restricts to an equivalence
$\cO_{\!X,y}\Mod_\coh\isom\cO_{\!X,y}\Mod^o_\coh$, and there
is a natural isomorphism
\set\begin{equation}\label{eq_noam-is-back}
R^i\Gamma_{\!\{y\}}C^\bullet_{|X(y)}\isom
D_y(H^{c(y)-i}(X(y),(\cD C^\bullet)_{|X(y)}))
\qquad
\text{in $\cO_{\!X,y}\Mod$}
\end{equation}
where $\cD$ is the duality functor corresponding to $\omega^\bullet$.
Notice that $a\cdot\one_{D_yM}=D_y(a\cdot\one_M)$ for every
$a\in\cO_{\!X,y}$, so \eqref{eq_noam-is-back} and
\eqref{eq_tough-noam} imply that (f) holds if and only if
\set\begin{equation}\label{eq_transform}
I\cdot H^{c-i-1}(X(y),(\cD C^\bullet)_{|X(y)})=0
\qquad
\text{for every closed point $y$ of $U$}.
\end{equation}
However, if $A$ is as in case (b) (resp. as in case (a))
of \eqref{subsec_loc-duality}, lemma \ref{lem_cod-function}(iii)
(resp. theorem \ref{th_dual-give-catenary}) says that every point
of $U$ specializes to a closed point, so \eqref{eq_transform}
holds if and only if $I\cdot H^{c-i-1}(\cD C^\bullet)_{|U}=0$.
Since $\cD C^\bullet$ is a complex of coherent $\cO_{\!X}$-module,
the latter condition is equivalent to saying that there exists
an ideal $I'$ whose radical is the maximal ideal of $A$, and
such that $I'I\cdot H^{c-i-1}(\cD C^\bullet)=0$. By invoking
\eqref{eq_noam-is-back} with $y:=x$, we see that this holds
if and only if $I'I\cdot R^{i+1}_{\{x\}}C^\bullet=0$. Summing up,
we have shown that (e)$\Leftrightarrow$(f).

\begin{claim}\label{cl_multiply-by-I}
The $A$-module $I\cdot R^i\Gamma_{\!\{x\}}C^\bullet$
is finitely copresented, for every $i\in\Z$.
\end{claim}
\begin{pfclaim} Notice first that, due to \eqref{eq_noam-is-back},
the $A$-module $E:=R^i\Gamma_{\!\{x\}}C^\bullet$ is finitely
copresented, so the same holds for $E^{\oplus n}$, for every
$n\in\N$. Since $I$ is finitely generated, it is the
image of a morphism $f:A^{\oplus n}\to A$, for some $n\in\N$,
and then $IE$ is naturally identified with the image of
$f\otimes_A\one_E$. Say that $E=D_x(M)$ for a finitely
presented $A$-module $M$; then clearly
$f\otimes\one_E=D_x(f^\vee\otimes_A\one_M)$, where
$f^\vee:A\to A^{\oplus n}$ is the transpose of $f$. Now
the assertion follows from claim \ref{cl_both-copresent}(i).
\end{pfclaim}

Taking into account claim \ref{cl_multiply-by-I} and lemma
\ref{lem_copresentation}, we deduce that
(c)$\Leftrightarrow$(d)$\Leftrightarrow$(e).

Obviously, (a)$\Rightarrow$(b). Next, consider the exact
sequence
$$
H^i(X,C^\bullet)\xrightarrow{\ \alpha_i\ }
H^i(U,C^\bullet)\xrightarrow{\ \beta_i\ }R^{i+1}\Gamma_{\!\{x\}}C^\bullet
\xrightarrow{\ \gamma_{i+1}\ }H^{i+1}(X,C^\bullet)
$$
(see \eqref{subsec_Gamma_Z}). Assume (b), which implies that
$I\cdot\Img\,\beta_i$ is a finitely generated $A$-module,
so there exists an ideal $I'\subset A$ whose radical is
the maximal ideal, and such that $I'I\cdot\Img\,\beta_i=0$.
On the other hand, by claim \ref{cl_multiply-by-I} and
proposition \ref{prop_copresentation} we know that
$\Img\,\gamma_{i+1}$ is a finitely presented $A$-module,
hence there exists an ideal $I''\subset A$ whose radical is
the maximal ideal, and such that $I''\cdot\Img\,\gamma_{i+1}=0$.
We conclude that $I''I'I\cdot R^{i+1}\Gamma_{\!\{x\}}C^\bullet=0$,
whence (b)$\Rightarrow$(e).

Lastly, we suppose that (e) holds, and we show that (a) follows.
To this aim, set $M:=I\cdot H^i(U,C^\bullet)\cap\Img\,\alpha_i$,
and consider the induced ladder with exact rows
\set\begin{equation}\label{eq_kyrgistan}
{\diagram
0 \ar[r] & M \ar[r] \ar[d] & I\cdot H^i(U,C^\bullet) \ar[r] \ar[d]
& I\cdot\Ker\,\gamma_{i+1} \ar[d] \ar[r] & 0 \\
0 \ar[r] & \Coker\,\gamma_i \ar[r]^-{\bar\alpha_i} &
H^i(U,C^\bullet) \ar[r] & \Ker\,\gamma_{i+1} \ar[r] & 0.
\enddiagram}
\end{equation}
From \eqref{eq_noam-is-back} we know that $\Ker\,\gamma_{i+1}$
is finitely copresented and arguing as in the proof of
claim \ref{cl_multiply-by-I}, we conclude that the same
holds for $I\cdot\Ker\,\gamma_{i+1}$. But (e) implies
that the latter $A$-module is annihilated by a finitely
generated ideal whose radical is the maximal ideal of $A$,
so it is finitely presented, by lemma \ref{lem_copresentation}.
We are thus reduced to checking that $M$ is a finitely
presented $A$-module. Now, set
$N:=\Ker\,(\bar\alpha_i\otimes_A\one_{A/I})$
and notice the short exact sequence
$$
0\to I\cdot\Coker\,\gamma_i\to M\to N\to 0.
$$
By proposition \ref{prop_copresentation}(ii) and
\eqref{eq_noam-is-back}, the $A$-module $\Coker\,\gamma_i$
is finitely presented, so the same holds for
$I\cdot\Coker\,\gamma_i$, because $A$ is coherent and
$I$ is finitely generated. So, we are further reduced
to showing that $N$ is finitely presented. However,
$N$ is naturally identified with the image of the map
$$
\Tor_1^A(\Ker\,\gamma_{i+1},A/I)\to(\Coker\,\gamma_i)\otimes_AA/I
$$
induced by the bottom row of \eqref{eq_kyrgistan}. On the
other hand, we have already noticed that $\Coker\,\gamma_i$
is finitely presented, so the same holds for
$(\Coker\,\gamma_i)\otimes_AA/I$; also, proposition
\ref{prop_copresentation}(i) and \eqref{eq_noam-is-back}
say that $\Ker\,\gamma_{i+1}$ is finitely copresented. Taking
into account proposition \ref{prop_copresentation}(ii) it
finally suffices to show :

\begin{claim} If $Q$ is any finitely copresented $A$-module
and $P$ any finitely presented $A$-module, the $A$-module
$\Tor_i^A(Q,P)$ is finitely copresented as well, for every
$i\in\N$.
\end{claim}
\begin{pfclaim}[] Since $I$ is finitely generated and
$A$ is coherent, we may find a resolution $F_\bullet\isom P[0]$
consisting of free $A$-modules of finite type (details
left to the reader), and then
$\Tor_i^A(Q,P)\simeq H_i(F_\bullet\otimes_AQ)$, so the
assertion follows easily from claim \ref{cl_both-copresent}(i)
(details left to the reader).
\end{pfclaim}
\end{proof}

\begin{proposition}\label{prop_flat-crit}
In the situation of \eqref{subsec_loc-duality}, let $d\in\N$
be any integer, and suppose additionally that :
\begin{itemize}
\item[(a')]
If $A$ is as in case {\em (a)}, then $A$ is regular of
dimension $d+1$.
\item[(b')]
If $A$ is as in case {\em (b)}, the structure morphism
$f:X\to S$ is flat, $f(x)=s$ and $f^{-1}(s)$ is a
regular scheme of dimension $d$.
\end{itemize}
Let $j:U:=X\!\setminus\!\{x\}\to X$ be the open immersion,
$\cF$ a flat quasi-coherent $\cO_{\!U}$-module.
Then the following two conditions are equivalent :
\begin{enumerate}
\alphaenu\addenu\addenu
\item
$\Gamma(U,\cF)$ is a flat $A$-module.
\item
$\delta(x,j_*\cF)\geq d+1$.
\end{enumerate}
\end{proposition}
\begin{proof} If $d=0$, then $A$ is a valuation ring
(proposition \ref{prop_integral-fp-ext}(ii)), and $U$ is
the spectrum of the field of fractions of $A$, in which case
both (c) and (d) hold trivially. Hence we may assume that $d\geq 1$.
Suppose now that (c) holds, and set $F:=\Gamma(U,\cF)$.
In view of \eqref{eq_Gamma-as-a-j}, we need to show that
$H^i(U,\cF)=0$ whenever $1\leq i\leq d-1$.
If $F^\sim$ denotes the $\cO_{\!X}$-module determined by $F$,
then $F^\sim_{|U}=\cF$. Moreover, by \cite[Ch.I, Th.1.2]{La}, $F$
is the colimit of a filtered family $(L_i~|~i\in I)$ of
free $A$-modules of finite rank, hence
$H^i(U,\cF)=\colim_{i\in I}H^i(U,L^\sim_i)$ by proposition
\ref{prop_dir-im-and-colim}(ii), so we are reduced to the
case where $\cF=\cO_{\!U}$, and therefore $j_*\cF=\cO_{\!X}$
(corollary \ref{cor_depth-cons}). In case (a'), assertion
(d) follows already, by virtue of \cite[Th.17.8]{Mat}.
So, suppose that (b') holds; since $f^{-1}(s)$ is regular, we have
$\delta(x,\cO_{\!f^{-1}(s)})=d$; then, since the topological space
underlying $X$ is noetherian, lemma \ref{lem_ineq-deltas}(ii)
and corollary \ref{cor_depth-flat-basechange} imply that
$\delta(x,\cO_{\!X})\geq d+1$, which is (d).

Conversely, suppose that (d) holds; we consider first the
case (b'), and we shall derive (c) by induction on $d$, the
case $d=0$ having already been dealt with. Let $\fm_A$ (resp.
$\bar\fm_A$) be the maximal ideal of $A$ (resp. of
$\bar A:=A/\fm_K A$). Suppose then that $d\geq 1$ and that
the assertion is known whenever $\dim f^{-1}(s)<d$. Pick
$\bar t\in\bar\fm_A\!\setminus\!\bar\fm_A^2$, and let $t\in\fm_A$
be any lifting of $\bar t$. Since $f^{-1}(s)$ is regular,
$\bar t$ is a regular element of $\bar A$, so that $t$ is
regular in $A$ and the induced morphism $g:X':=\Spec\,A/tA\to S$
is flat (\cite[Ch.IV, Th.11.3.8]{EGAIV-3}). Let $j':U':=U\cap X'\to X'$
be the restriction of $j$; our choice of $\bar t$ ensures that
$g^{-1}(s)$ is a regular scheme, so the pair $(X',\cF/t\cF_{|U'})$
fulfills the conditions of the proposition, and $\dim g^{-1}(s)=d-1$.
Furthermore, since $\cF$ is a flat $\cO_{\!U}$-module, the sequence
$0\to\cF\xrightarrow{t}\cF\to\cF/t\cF\to 0$ is short exact.
Assumption (d) means that
\set\begin{equation}\label{eq_H^i-vanish}
H^i(U,\cF)=0\qquad\text{whenever $1\leq i\leq d-1$}.
\end{equation}
Therefore, from the long exact cohomology sequence we deduce
that $H^i(U,\cF/t\cF)=0$ for $1\leq i\leq d-2$, {\em i.e.}
\set\begin{equation}\label{eq_about-delta}
\delta(x,j'_*\cF/t\cF)\geq d.
\end{equation}
The same sequence also yields a left exact sequence :
\set\begin{equation}\label{eq_include-mod-t}
0\to H^0(U,\cF)\otimes_AA/tA\xrightarrow{\alpha}
H^0(U,\cF/t\cF)\to H^1(U,\cF).
\end{equation}

\begin{claim}\label{cl_flat-over-t}
$H^0(U,\cF)\otimes_AA/tA$ is a flat $A/tA$-module.
\end{claim}
\begin{pfclaim}
By \eqref{eq_about-delta} and  our inductive assumption,
$H^0(U,\cF/t\cF)$ is a flat $A/tA$-module. In case $d=1$,
proposition \ref{prop_integral-fp-ext}(ii) shows that $A/tA$
is a valuation ring, and then the claim follows from
\eqref{eq_include-mod-t}, and \cite[Ch.VI, \S3, n.6, Lemma 1]{BouAC}.
If $d>1$, \eqref{eq_H^i-vanish} implies that $\alpha$
is an isomorphism, so the claim holds also in such case.
\end{pfclaim}

Set $V:=\Spec\,A[t^{-1}]=X\!\setminus\! V(t)\subset U$;
since $\cF_{|V}$ is a flat $\cO_V$-module, the
$A[t^{-1}]$-module
$H^0(U,\cF)\otimes_AA[t^{-1}]\simeq H^0(V,\cF_{|V})$
is flat. Moreover, since $t$ is regular on both $A$ and
$H^0(U,\cF)$, an easy calculation shows that
$\Tor_i^A(A/tA,H^0(U,\cF))=0$ for $i>0$. Then the
contention follows from claim \ref{cl_flat-over-t}
and \cite[Lemma 5.2.1]{Ga-Ra}.

Lastly, in case (a') holds, one picks any element
$t\in\fm_A\!\setminus\!\fm^2_A$ and argues in the same way
(with some simplifications : details left to the reader).
\end{proof}

\begin{remark} The special case of proposition
\ref{prop_flat-crit} where $A$ is regular local ring and
$\cF$ is a locally free $\cO_{\!U}$-module of finite
rank has been studied in detail in \cite{Hor}.
\end{remark}

\begin{lemma}\label{lem_easy-short}
Let $A$ be a noetherian ring, $M$ an $A$-module
endowed with a finite filtration
$M_0:=0\subset M_1\subset\cdots\subset M_k:=M$,
and denote by $\gr^\bullet M$ the associated graded
$A$-module. Suppose that, for every $i=1,\dots,k$,
there exists an ideal $I_i\subset A$ such that
$I_i\cdot\gr^iM$ is an $A$-module of finite type.
Then $\prod_{i=1}^kI_iM$ is an $A$-module of finite
type.
\end{lemma}
\begin{proof} A simple induction on the length of the
filtration reduces the assertion to the following 

\begin{claim}\label{cl_easy-short}
Let $A$ be a noetherian ring, $I,J\subset A$ two ideals, and
$$
0\to M_1\to M_2\to M_3\to 0
$$
a short exact sequence of $A$-modules, such that $IM_1$
and $JM_3$ are finitely generated. Then $IJM_2$ is finitely
generated.
\end{claim}
\begin{pfclaim}[] Indeed, the induced map
$N:=JM_2/(M_1\cap JM_2)\to JM_3$ is injective, so $N$ is
finitely generated, and then the same holds for $I\otimes_AN$;
but the latter maps surjectively onto
$IJM_2/(IM_1\cap IJM_2)\simeq(IJM_2+IM_1)/IM_1$, so the
latter is also of finite type, and therefore the same holds
for $IJM_2+IM_1$, whence the lemma.
\end{pfclaim}
\end{proof}

\sset\subsubsection{}\label{subsec_SGA-2-fincrit}
Let $A$ be a noetherian ring, $\omega^\bullet$ a dualizing
complex on $X:=\Spec\,A$, and $I\subset A$ any ideal.
Let also $U\subset X$ be any open subset, and $K^\bullet$
any object of $\sD^b(\cO_{\!U}\Mod)_\coh$. We denote :
\begin{itemize}
\item
$\bar U$ the topological closure of $U$ in $X$
\item
$Z:=X\!\setminus\!U$ and $\partial U:=\bar U\setminus U$,
endowed with the reduced closed subscheme structures
\item
$U_1:=\{x\in U~|~\text{there exists $y\in\partial U$ which
is an immediate specialization of $x$ in $\bar U$}\}$
\item
$\cD_U:\sD^b(\cO_{\!U}\Mod)_\coh^o\to\sD^b(\cO_{\!U}\Mod)_\coh$
the duality functor associated with $\omega^\bullet_{|U}$
\item
$n(K^\bullet)$ the cardinality of
$\{i\in\Z~|~H^i(\cD_U K^\bullet)\neq 0\}$
\item
$d(K^\bullet):=\min(n(K^\bullet),1+\dim\partial U)$
\item
$\partial U(x):=U\cap X(x)$ and $\delta(x):=\dim\partial U(x)$
for every $x\in X$
\item
$e(x):=\min(1+\delta(y)~|~
\text{$y$ is a specialization of $x$ in $\bar U$})$ for every
$x\in U_1$.
\end{itemize}
Here we use the convention that $\dim\partial U=-1$ if
$\partial U=\emptyset$, and likewise we define the dimension
of $\partial U(y)$, for any $y\in X$ such that this set is
empty. On the other hand, the dimension of $\partial U$ can
be infinite, in the current setting. Notice that the integer
$n(K^\bullet)$ depends only on $K^\bullet$, and does not depend
on the choice of $\omega^\bullet$.

\begin{lemma}\label{lem_closed-points}
For every $x\in\partial U$, every closed point of\/ $U\cap X(x)$
is an immediate generization of $x$ in $X(x)$. Especially,
every such point lies in $U_1$.
\end{lemma}
\begin{proof} Fix any $x\in\partial U$, set $A:=\cO_{\!X,x}$,
and for any closed point $z$ of $U\cap X(x)$, denote by
$\fp_z\subset A$ the prime ideal corresponding to $z$, so
that $\fp_z$ is strictly contained in the maximal ideal
$\fm_A$ of $A$. Let also $I\subset A$ be any ideal such
that $\Spec\,A/I=Z(x)$, so $I\not\subset\fp_z$ by construction.
Set $B:=A/\fp_z$, and denote by $\fm_B,\bar I\subset B$
the images of respectively $\fm_A$ and $I$. It suffices
to show that $\dim B=1$; suppose by contradiction that
the latter fails, and pick any $f\in\bar I\setminus\{0\}$.

\begin{claim}\label{cl_Ratliff}
If $\dim B>1$, there exists a prime ideal $\bar\fq\subset B$
of height $1$, such that $f\notin\bar\fq$.
\end{claim}
\begin{pfclaim} Clearly, every prime ideal of height one
containing $f$ is a maximal point of $\Supp\,B/fB$, hence
it lies in $\Ass\,B/fB$ (lemma \ref{lem_ass-are-mins}(i));
therefore, the set $\Sigma$ of prime ideals of $B$ of height
one containing $f$ is finite (\cite[Th.6.5(i)]{Mat}). If
$\Sigma$ is empty, Krull's {\em Hauptidealsatz}
(\cite[Th.13.5]{Mat}) implies that $f$ must be invertible,
in which case the claim is obvious. If $\Sigma$ contains
a single prime ideal $\bar\fq_1$, notice that
$\bar\fq_1\neq\fm_B$, since $\dim B>1$, and pick any
$g\in\fm_B\setminus\fq_1$; again by Krull's theorem we
know that there exists a prime ideal $\bar\fq\subset B$
of height one containing $g$, and by construction
$\bar\fq\neq\bar\fq_1$, whence the claim, in this case.
Otherwise, say that $\Sigma=\{\bar\fq_1,\cdots,\bar\fq_n\}$
for some $n\geq 2$. For every $i,j\leq n$ with $i\neq j$ we
may find an element $g_{ij}\in\fq_i\setminus\fq_j$, and we
set $g:=\sum_{j=1}^n\prod_{i\neq j}g_{ij}$. It is easily seen
that $g\in\fm_B\setminus(\bar\fq_1\cup\cdots\cup\bar\fq_n)$,
so we may apply again Krull's theorem to exhibit a prime ideal
$\bar\fq$ which cannot lie in $\Sigma$, since it contains $g$.
\end{pfclaim}

Let $\bar\fq$ be as in claim \ref{cl_Ratliff}, and denote
by $\fq\subset A$ the preimage of $\bar\fq$; then $\fq$
corresponds to a point $y\in U\cap X(x)$ which is different
from $x$, and is a proper specialization of $z$ in $X(x)$,
which is impossible, since by assumption $z$ is closed in $U$.
\end{proof}

\begin{theorem}\label{th_SGA-2-fincrit}
In the situation of \eqref{subsec_SGA-2-fincrit},
let $i\in\Z$ be any integer, and consider the following
conditions :
\begin{enumerate}
\alphaenu
\item
$I\cdot R^i\Gamma_{\!\{x\}}K^\bullet_{|U(x)}=0$ for every $x\in U_1$.
\item
$I^{d(K^\bullet)}\cdot H^i(U,K^\bullet)$ is an $A$-module of
finite type.
\item
$I\cdot H^j(U,K^\bullet)$ is an $A$-module of finite type,
for every $j\leq i$.
\item
$I^{e(x)}\cdot R^j\Gamma_{\!\{x\}}K^\bullet_{|U(x)}=0$ for
every $x\in U_1$ and every $j\leq i$.
\end{enumerate}
Then {\em(a)}$\Rightarrow${\em(b)} and {\em(c)}$\Rightarrow${\em(d)}.
\end{theorem}
\begin{proof}(a)$\Rightarrow$(b): By corollary \ref{cor_fp-approx},
we may find an object $L^\bullet$ of $\sD^b(\cO_{\!X}\Mod)_\coh$ with
an isomorphism $L^\bullet_{|U}\isom K^\bullet$ in
$\sD^b(\cO_{\!U}\Mod)_\coh$, whence an exact sequence
$$
H^i(X,L^\bullet)\to H^i(U,K^\bullet)\to R^{i+1}\Gamma_{\!Z}L^\bullet
\to H^{i+1}(X,L^\bullet)
$$
(see \eqref{subsec_Gamma_Z}) whose first and last terms are
finitely generated $A$-modules. By an easy diagram chase, it
follows that condition (b) is equivalent to
\begin{enumerate}
\addenu\addenu\addenu\addenu
\alphaenu
\item
$I^{d(K^\bullet)}\cdot R^{i+1}\underline\Gamma_{\!Z}L^\bullet$ is
a coherent $\cO_{\!X}$-module.
\end{enumerate}
It will therefore suffice to show that (a)$\Rightarrow$(e).
However, from lemma \ref{lem_adj-Gamma_Z}(iii) we get a natural
isomorphism
$$
R\underline\Gamma_{\!Z}L^\bullet\isom
R\cHom_Z^\bullet(\cD_X L^\bullet,\omega^\bullet)
$$
where $\cD_X$ is the dualizing functor on $\sD^b(\cO_{\!X}\Mod)_\coh$
defined by $\omega^\bullet$.
By choosing an injective resolution for $\omega^\bullet$,
we obtain a spectral sequence
$$
E_2^{pq}:=R^p\cHom^\bullet_Z(H^{-q}(\cD_X L^\bullet)[0],\omega^\bullet)
\Rightarrow R^{p+q}\cHom_Z^\bullet(\cD_X L^\bullet,\omega^\bullet).
$$
Now, for every $q\in\Z$, define the $\cO_{\!X}$-module
$\bar H{}^{-q}$ by the short exact sequence
$$
0\to\underline\Gamma_ZH^{-q}(\cD_XL^\bullet)\to
H^{-q}(\cD_XL^\bullet)\to\bar H{}^{-q}\to 0.
$$
There results an exact sequence
\set\begin{equation}\label{eq_EFG}
F^{pq}:=R^p\cHom^\bullet_Z(\bar H{}^{-q}[0],\omega^\bullet)\to
E^{pq}_2\to
G^{pq}:=
R^p\cHom^\bullet_Z
(\underline\Gamma_ZH^{-q}(\cD_XL^\bullet)[0],\omega^\bullet)
\end{equation}

\begin{claim}\label{cl_finite-piece}
For every $p,q\in\Z$ we have :
\begin{enumerate}
\item
$E^{pq}_2$ and $F^{pq}$ are quasi-coherent $\cO_{\!X}$-modules.
\item
$G^{pq}$ and $\bar H{}^{-q}$ are coherent $\cO_{\!X}$-modules.
\end{enumerate} 
\end{claim}
\begin{pfclaim} The assertions about $E^{pq}_2$ and $F^{pq}$ follow
immediately from lemmata \ref{lem_without-cohereur}(i) and
\ref{lem_adj-Gamma_Z}(iii) and corollary \ref{cor_Ext-loc=glob}(ii),
and the assertion for $\bar H{}^{-q}$ follows from lemma
\ref{lem_without-cohereur}(i). Next, notice that the natural map
$$
\underline\Gamma_ZH^{-q}(\cD_XL^\bullet)\to R\underline\Gamma_Z
\circ\underline\Gamma_{\!Z}H^{-q}(\cD_XL^\bullet)
$$
is an isomorphism in $\sD(\cO_{\!X}\Mod)$, whence an isomorphism
of $\cO_{\!X}$-modules
\set\begin{equation}\label{eq_target-this}
G^{pq}\isom R^p\cHom^\bullet_{\cO_{\!X}}
(\underline\Gamma_ZH^{-q}(\cD_XL^\bullet)[0],\omega^\bullet)
\end{equation}
(lemma \ref{lem_adj-Gamma_Z}(vi)). But since $X$ is noetherian,
$\underline\Gamma_ZH^{-q}(\cD_XL^\bullet)$ is a coherent
$\cO_{\!X}$-module, so the same holds for the target of
\eqref{eq_target-this} (proposition \ref{prop_replace-prop.12.3.5}(i)),
and the claim follows.
\end{pfclaim}

\begin{claim}\label{cl_to-prove-e}
In order to show (e), it suffices to check that :
\begin{enumerate}
\item
$I\cdot F^{pq}=0$ for every $p,q\in\Z$ such that $p+q=i+1$.
\item
$F^{pq}_x=0$ for every $x\in X$ and every $p\in\Z$
such that $p\notin[c_X(x)-\dim\partial U(x),c_X(x)]$.
\end{enumerate}
\end{claim}
\begin{pfclaim}
Say that $\cD_XL^\bullet\in\sD^{[a,b]}(\cO_{\!X}\Mod)$
for some $a,b\in\Z$ with $a\leq b$; to begin with, we show that
-- under assumptions (i) and (ii) of the claim -- there exist
for every $q=a,\dots,b$ and every
$p=i+1-b,\dots,i+1-a$ and $x\in X$, integers
$\nu(q),\mu(p,x)\in\N$ such that $I^{\nu(q)}F^{pq}=0$
and $I^{\mu(p,x)}F^{pq}_x=0$, and moreover
\set\begin{equation}\label{eq_nu-mu-bounds}
\sum_{q=a}^b\nu(q)\leq n(K^\bullet)
\qquad
\sum_{p=i+1-b}^{i+1-a}\mu(p,x)\leq 1+\delta(x)
\quad
\text{for every $x\in X$}.
\end{equation}
Indeed, if $q\in[a,b]$ is an index such that
$H^{-q}(\cD_UK^\bullet)=0$, it follows that
$H^{-q}(\cD_XL^\bullet)_{|U}=0$, so $\bar H{}^{-q}=0$,
therefore $F^{pq}=0$ for every $p\in\Z$, and we may take
$\nu(q)=0$ in this case. But assumption (i) says that
in any case we may take $\nu(q)\leq 1$ for every
$q=a,\dots,b$. This already suffices to achieve
the sought upper bound for the sum of the values $\nu(q)$.
The other upper bound is treated similarly : we may in any
case take $\mu(p,x)\leq 1$, due to condition (i), and moreover
condition (ii) says that we may take $\mu(p,x)=0$ for every
$p\in\Z$ such that $p\notin[c_X(x)-\dim\partial U(x),c_X(x)]$,
whence the contention.

Next, recall that $M:=R^{i+1}\underline\Gamma_ZL^\bullet$ admits
a finite filtration $\Fil^\bullet M$ whose associated graded
$\cO_{\!X}$-module is $\bigoplus_{p+q=i+1}E^{pq}_\infty$, and by
construction, each term $E^{pq}_\infty$ is a quasi-coherent
subquotient of $E^{pq}_2$ (claim \ref{cl_finite-piece}(i)),
so $\Fil^qM$ is a quasi-coherent $\cO_{\!X}$-module as well,
for every $q=a,\dots,b$.
Using claims \ref{cl_easy-short} and \ref{cl_finite-piece}(i,ii)
we deduce that $I^{\nu(q)}E^{pq}_2$ is a coherent
$\cO_{\!X}$-module for every $q=a,\dots,b$, so the same
holds for $I^{\nu(q)}E^{pq}_\infty$, and then
\eqref{eq_nu-mu-bounds} and lemma \ref{lem_easy-short}
imply that $I^{n(K^\bullet)}M$ is a coherent $\cO_{\!X}$-module.

Set $\delta:=\dim\partial U\in\N\cup\{\infty\}$; if
$\delta\geq n(K^\bullet)-1$, the proof is complete; otherwise
$\delta\in\N$, and it remains only to check that $I^{1+\delta}M$
is also a coherent $\cO_{\!X}$-module. To this aim, notice
that \eqref{eq_EFG} induces a short exact sequence
$$
0\to F_\infty^{pq}\to E^{pq}_\infty\to G^{pq}_\infty\to 0
$$
where $F^{pq}_\infty$ (resp. $G^{pq}_\infty$) is a certain
quasi-coherent subquotient of $F^{pq}$ (resp. of $G^{pq}$);
from claim \ref{cl_finite-piece}(ii) it follows that
$G^{pq}_\infty$ is a coherent $\cO_{\!X}$-module, and the
foregoing implies that $I^{\mu(p,x)}F^{pq}_{\infty,x}=0$ for
every $x\in X$ and every $p=i+1-b,\dots,i+1-a$. We may
then find, for every $q=a,\dots,b$ a coherent
$\cO_{\!X}$-submodule $N^q\subset\Fil^qM$ such that
the induced morphism $\Fil^qM\to G^{pq}_\infty$ restricts
to an epimorphism $N^q\to G^{pq}_\infty$. Set
$N:=\sum_{q=a}^bN^q\subset M$, and $P:=M/N$, endow
$P$ with the filtration $\Fil^\bullet P$ induced by the
filtration of $M$, and denote by $\gr^\bullet P$ the
associated graded $\cO_{\!X}$-module. By construction,
we get an epimorphism $\psi_q:E^{pq}_\infty\to\gr^qP$ for
every $q=a,\dots,b$, whose kernel contains the image
of $N^q$ in $E^{pq}_\infty$. Consequently, the restriction
of $\psi_q$ to $F^{pq}_\infty$ is still an epimorphism,
so that $I^{\mu(p,x)}\gr^{i+1-p}P_x=0$ for every $x\in X$
and every $p=i+1-b,\dots,i+1-a$. In view of
\eqref{eq_nu-mu-bounds} we conclude that
$I^{1+\delta(x)}P_x=0$ for every $x\in X$, and since clearly
$\delta(x)\leq\delta$ for every $x\in X$, it follows that
$I^{1+\delta}P=0$. Lastly, since $N$ is a coherent
$\cO_{\!X}$-module, it suffices to invoke claim
\ref{cl_easy-short} to conclude the proof.
\end{pfclaim}

Let $\cR^\bullet$ be the residual complex arising from
$\omega^\bullet$, introduced in \eqref{subsec_residual}.
With this notation, $F^{pq}$ is a subquotient of the
$A$-module
$\cHom_{\cO_{\!X}}(\bar H{}^{-q},\underline\Gamma_Z\cR^p)$.
From lemma \ref{lem_flabby-Gamma_Z}(iii.b), we get
$$
\underline\Gamma_Z\cR^p=\bigoplus_{x\in X_p\cap Z}E(x)^\sim.
\qquad
\text{with $X_p:=\{x\in X~|~c_X(x)=p\}$
for every $p\in\Z$}
$$
where $c_X:|X|\to\Z$ is the weak codimension function associated
with $\omega^\bullet$ as in \eqref{subsec_dual-catenary} and
$E(x)$ is an injective hull of the $A$-module $\kappa(x)$,
for every $x\in X$ (notation of \eqref{sec_various-O-mod}).
Hence, denote by $X^*_p$ the set of all specializations in
$X$ of the points of $X_p$; we deduce
\set\begin{equation}\label{eq_precise-support}
\Supp\,\underline\Gamma_Z\cR^p=Z\cap X^*_p
\end{equation}
(we caution the reader that $X_p$ is usually not a
pro-constructible subset of $X$, so $X^*_p$ is usually
strictly contained in the topological closure of $X_p$
in $X$). Moreover, since $\bar H{}^{-q}$ is a coherent
$\cO_{\!X}$-module (claim \ref{cl_finite-piece}(ii))
and $\underline\Gamma_Z\cR^p$ is quasi-coherent
(lemma \ref{lem_without-cohereur}(i)), we also see that
\set\begin{equation}\label{eq_support-homs}
\cHom_{\cO_{\!X}}(\bar H{}^{-q},\underline\Gamma_Z\cR^p)_x\simeq
\Hom_{\cO_{\!X,x}}(\bar H{}^{-q}_x,(\underline\Gamma_Z\cR^p)_x)
\qquad
\text{for every $x\in X$}
\end{equation}
(\cite[Lemma 2.4.29(i.a)]{Ga-Ra}). Combining
\eqref{eq_precise-support} and \eqref{eq_support-homs}
we conclude that
$$
\Supp\,F^{pq}\subset Z\cap X^*_p\cap\Supp\,\bar H{}^{-q}.
$$
From this, we can already show condition (ii) of claim
\ref{cl_to-prove-e} : indeed, say that $F^{pq}_x\neq 0$;
then $x\in X^*_p$, so $c_X(x)\geq p$, and on the other
hand, by construction we have
$\underline\Gamma_Z\bar H{}^{-q}=0$, so $\Supp\,H{}^{-q}$
is the topological closure in $X$ of
$U\cap\Supp\,\bar H{}^{-q}$, therefore $x\in\partial U$
and furthermore there exists a generization $z$ of $x$
in $X$ with $c_X(z)=p$, so $z\in X_p\cap\partial U(x)$,
which forces the bound $c_X(x)\leq p+\dim\partial U(x)$,
as required. 

Lastly, to get condition (i) of claim
\ref{cl_to-prove-e} it suffices to check assertions (i)
and (ii) of the following :

\begin{claim}
(i)\ \ $\Supp\,(I\cdot F^{pq})\subset
\Supp\,\cHom_{\cO_{\!X}}(I\cdot\bar H{}^{-q},\underline\Gamma_Z\cR^p)$
for every $p,q\in\Z$.
\begin{enumerate}
\addenu
\item
$\cHom_{\cO_{\!X}}(I\cdot\bar H{}^{-q},\underline\Gamma_Z\cR^p)=0$
whenever $p+q=i+1$.
\item
$Z\cap X_p\cap\Supp\,(I\cdot\bar H{}^{-q})=\emptyset$
whenever $p+q=i+1$.
\end{enumerate}
\end{claim}
\begin{pfclaim}(i): Consider arbitrary $x\in X$, $a\in I$ and
$f\in F^{pq}_x$; by virtue of \eqref{eq_support-homs}, the
element $f$ is represented by an $\cO_{\!X,x}$-linear map 
$\phi:\bar H{}^{-q}_x\to\underline\Gamma_Z\cR^p_x$, and therefore
$af$ is represented by $a\phi$; the latter obviously factors
through the submodule $I\cdot \bar H{}^{-q}_x$, so the
contention follows from the natural identification
\set\begin{equation}\label{eq_I-support-homs}
\cHom_{\cO_{\!X}}(I\cdot\bar H{}^{-q},\underline\Gamma_Z\cR^p)_x\simeq
\Hom_{\cO_{\!X,x}}(I\cdot\bar H{}^{-q}_x,(\underline\Gamma_Z\cR^p)_x)
\qquad
\text{for every $x\in X$}
\end{equation}
which is proven as \eqref{eq_support-homs}.

(ii): Let $x\in X$ be any point, and
$\phi:I\cdot\bar H{}^{-q}_x\to(\underline\Gamma_Z\cR^p)_x$
any $\cO_{\!X,x}$-linear map; taking into account
\eqref{eq_I-support-homs}, we are reduced to showing that
$\phi=0$. Now, notice that
$$
(\underline\Gamma_Z\cR^p)_x=
\bigoplus_{y\in Z(x)_p}E(y)
\qquad
\text{where $Z(x)_p:=X_p\cap Z\cap X(x)$}.
$$
If $Z(x)_p=\emptyset$, we are done; otherwise, for every
$y\in Z(x)_p$ denote by $p_y:(\underline\Gamma_Z\cR^p)_x\to E(y)$
the natural projection; it suffices to check that
$p_y\circ\phi=0$ for every such $y$. However,
$p_y\circ\phi$ factors through an $\cO_{\!X,y}$-linear
map $I\cdot\bar H{}^{-q}_y\to(\underline\Gamma_Z\cR^p)_y$,
so we may replace $x$ by $y$, and assume from start that
$x\in X_p\cap Z$, in which case the contention will
follow from assertion (iii) of the claim.

(iii): Fix any point $y\in X_p\cap Z$. We have already
noticed that $\Supp\,\bar H{}^{-q}$ is the topological
closure in $X$ of $U\cap\Supp\,\bar H{}^{-q}$; we are
then reduced to showing that $I\cdot\bar H{}^{-q}_z=0$
for every closed point $z$ of $X(y)\cap U$. But for any
such $z$, we have an isomorphism
\set\begin{equation}\label{eq_dai-dai}
H^{-q}(\cD_UK^\bullet)_z\isom\bar H{}^{-q}_z
\end{equation}
of $\cO_{\!X,z}$-modules. On the other hand, corollary
\ref{cor_loc-duality}(iii) yields an isomorphism of
$\cO_{\!X,z}$-modules
\set\begin{equation}\label{eq_uff-tired}
\cO_{\!X,z}^\wedge\otimes_{\cO_{\!X,z}}H^{-q}(\cD_UK^\bullet)_z
\isom D_z\circ R^{q+c_X(z)}\Gamma_{\!\{z\}}K^\bullet_{|U(z)}
\end{equation}
where $\cO_{\!X,z}^\wedge$ denotes the completion of the
ring $\cO_{\!X,z}$ and $D_z$ denotes the local duality
functor corresponding to the point $z$. We remark now that
\set\begin{equation}\label{eq_immediate-gener}
c_X(z)=c_X(y)-1=p-1.
\end{equation}
Indeed, in view of lemma \ref{lem_dual-catenary}, identity
\eqref{eq_immediate-gener} asserts that $y$ is an immediate
specialization of $z$ in $X(y)$, and this is already known
from lemma \ref{lem_closed-points}.

Finally, let $a\in I$ be any element; taking into account
\eqref{eq_dai-dai}, we are reduced to checking that scalar
multiplication by $a$ is the zero endomorphism
of $H^{-q}(\cD_UK^\bullet)_z$, which -- by virtue of
\eqref{eq_uff-tired} and \eqref{eq_immediate-gener} --
is equivalent to showing that scalar multiplication by $a$
is the zero endomorphism on 
$R^i\Gamma_{\!\{z\}}K^\bullet_{|U(z)}$, and the latter
is ensured by our condition (a).
\end{pfclaim}

(c)$\Rightarrow$(d): Let $x\in U_1$ be any point, pick
any immediate specialization $y\in\partial U$ of $x$ in
$X$ such that $1+\dim\partial U(y)=e(x)$, let
$U(y):=U\cap X(y)$, and notice that the natural map
$$
H^j(U,K^\bullet)\otimes_A\cO_{X,y}\to H^j(U(y),K^\bullet_{|U(y)})
$$
is an isomorphism for every $j\in\Z$, so we may replace
$X$ by $X(y)$, $U$ by $U(y)$ and $K^\bullet$ by
$K^\bullet_{|U(y)}$, and assume from start that $A$
is a local ring, and $y$ is the closed point of $X$.
In this situation, we need to show that
$I^eR^j\Gamma_{\!\{x\}}K^\bullet_{|U(x)}=0$, with
$e:=1+\dim\partial U$ and every $j\leq i$. Then,
let $V:=X\setminus\{y\}$, denote $j:U\to V$ the
induced open immersion, and set
$$
L^\bullet:=\tau^{\leq i}Rj_*K^\bullet
$$
(where $\tau^{\leq i}$ denotes the usual truncation functor),
so that $L\in\Ob(\sD^b(\cO_V\Mod)_\qcoh)$ and
$L^\bullet_{|U}=\tau^{\leq i}K^\bullet$. Notice also
that the natural morphism of $\cO_{\!X}$-modules
$$
H^j(U,K^\bullet)^\sim\to H^jL^\bullet
$$
is an isomorphism for every $j\leq i$. Taking into account
proposition \ref{prop_coher-then-noether}(ii), it follows
that we may find a morphism $\phi^\bullet:C^\bullet\to L^\bullet$
in $\sD^b(\cO_{\!X}\Mod)$ such that
\begin{itemize}
\item
$C^\bullet\in\Ob(\sD^b(\cO_{\!X}\Mod)_\coh)$.
\item
The induced morphism
$C^\bullet_{|U}\to\tau^{\leq i}K^\bullet_{|U}$ is an isomorphism
in $\sD^b(\cO_{\!U}\Mod)_\coh$.
\item
$H^j\phi^\bullet:H^jC^\bullet\to H^jL^\bullet$ is a monomorphism
for every $j\in\Z$, and its image contains the $\cO_{\!X}$-submodule
$I\cdot H^j(U,K^\bullet)^\sim$, for every $j\leq i$. 
\end{itemize}
Set $D^\bullet:=\Cone\,\phi$; by construction, for every $j\in\Z$
we have
\set\begin{equation}\label{eq_supports-L-D}
\Supp\,H^jL^\bullet\subset\bar U\cap V
\quad
\Supp\,H^jD^\bullet\subset Z\cap\Supp\,H^iL^\bullet\subset
\partial U\cap V
\quad
I\!\cdot\!H^jD^\bullet=0.
\end{equation}
Now, consider the spectral sequence
$$
E^{pq}_2:=H^p(V,H^qD^\bullet)\Rightarrow H^{p+q}(V,D^\bullet).
$$
It follows from \eqref{eq_supports-L-D} that $E^{pq}_2=0$ for
every $p>\dim\partial U\cap V=e-2$, and $I\cdot E^{pq}_2=0$
for every $p,q\in\Z$, whence
$$
I^{e-1}H^j(V,D^\bullet)=0
\qquad
\text{for every $j\in\Z$}.
$$
On the other hand, we have an exact sequence of $A$-modules
$$
H^{j-1}(V,D^\bullet)\to H^j(V,C^\bullet)\to H^j(V,L^\bullet)=
H^j(U,K^\bullet)
\qquad
\text{for every $j\leq i$}.
$$
In view of lemma \ref{lem_easy-short}, we conclude that
$I^eH^j(V,C^\bullet)$ is an $A$-module of finite type, whence
$I^e\cdot R^j\Gamma_{\!\{x\}}C^\bullet_{|U(x)}=0$ for every $j\leq i$
(theorem \ref{th_finit}). To conclude the proof, it now
suffices to observe that the natural map
$R^j\Gamma_{\!\{x\}}C^\bullet_{|U(x)}\to
R^j\Gamma_{\!\{x\}}K^\bullet_{|U(x)}$ is an isomorphism,
for every $j\leq i$.
\end{proof}

\begin{corollary} In the situation of \eqref{subsec_SGA-2-fincrit},
for every $i\in\Z$ we have :
\begin{enumerate}
\item
If $R^i\Gamma_{\!\{x\}}K^\bullet=0$ for every $x\in U_1$, then
$H^i(U,K^\bullet)$ is an $A$-module of finite type.
\item
The following conditions are equivalent :
\begin{enumerate}
\item
$R^j\Gamma_{\!\{x\}}K^\bullet=0$ for every $x\in U_1$
and every $j\leq i$.
\item
$H^j(U,K^\bullet)$ is an $A$-module of finite type for
every $j\leq i$.
\end{enumerate}
\end{enumerate}
\end{corollary}
\begin{proof} It is the special case of theorem
\ref{th_SGA-2-fincrit} with $I:=A$.
\end{proof}

\subsection{Hochster's theorem and Stanley's theorem}
\label{sec_homol-polytope}
The two main results of this section are theorems
\ref{th_Hochster} and \ref{th_Stanley}, which were proved
originally respectively by Hochster in \cite{Hoch}, and
by Stanley in \cite{Sta}. Our proofs follow in the main the
posterior methods presented by Bruns and Herzog in \cite{Br-He}
(which in turns, partially rely on some ideas of Danilov).
These results shall be used in section \ref{sec_log-regular},
in order to show that regular log schemes are Cohen-Macaulay.
We begin with some preliminaries from algebraic topology,
which we develop only to the extent that is required for
our special purposes: a much more general theory exists, and
is well known to experts (see {\em e.g.} \cite[Ch.IX]{Mas}).

For any topological space $X$, and any subset $T\subset X$,
we denote by $\bar T$ the topological closure of $T$ in $X$,
endowed with the topology induced from $X$. We let $H_\bullet(X)$
(resp. $H_\bullet(X,T)$) be the singular homology groups of $X$
(resp. of the pair $(X,T)$). Also, for every $n\in\N$ fix a
Banach norm $\Vert\cdot\Vert$ on $\R^n$, and for every real number
$\rho>0$, let $\B^n(\rho):=\{v\in\R^n~|~\Vert v\Vert<\rho\}$, and set
$\SS^{n-1}:=\bar\B{}^n(1)\setminus\B{}^n(1)$ (so $\SS^{-1}=\emptyset$).

\begin{definition}\label{def_cell-complex}
(i)\ \
A {\em finite regular cell complex} (or briefly : a
{\em cell complex}) is the datum of a topological
space $X$, together with a finite filtration
$$
X^{-1}:=\emptyset\subset X^0\subset X^1\subset X^2
\subset\cdots\subset X^k:=X
$$
consisting of closed subspaces, such that $X^0$ is a discrete
topological space, and for every $i=0,\dots,k$ we have a
decomposition
$$
X^i\setminus X^{i-1}=\bigcup_{\lambda\in\Lambda_i}e^i_\lambda
$$
where :
\begin{enumerate}
\alphaenu
\item
$\Lambda_i$ is a finite set, and for every $\lambda\in\Lambda_i$
there exists a homeomorphism
$$
f_\lambda:\bar\B{}^i(1)\isom\bar e{}^i_\lambda
\qquad
\text{such that $f^{-1}_\lambda(e^i_\lambda)=\B^i(1)$}.
$$
\item
$e^i_\lambda\cap e^i_\mu=\emptyset$ for any two distinct indices
$\lambda,\mu\in\Lambda_i$.
\item
$\bar e^i_\lambda\setminus e^i_\lambda\subset X^{i-1}$ for every
$\lambda\in\Lambda_i$.
\end{enumerate}

(ii)\ \
The smallest $k\in\N$ such that $X^k=X$ is called the
{\em dimension} of $X^\bullet$, and is denoted $\dim X^\bullet$.
For every $i=0,\dots,k$, the subsets $e^i_\lambda$ are called
the {\em $i$-dimensional cells of $X^\bullet$}.

(iii)\ \
A {\em subcomplex} of the cell complex $X^\bullet$ is a closed subset
$Y\subset X$ which is a union of cells of $X$. Then the filtration
such that $Y^i:=Y\cap X^i$ for every $i\geq -1$ defines a cell
complex structure $Y^\bullet$ on $Y$.
\end{definition}

\begin{lemma}\label{lem_cell-complex}
With the notation of definition {\em\ref{def_cell-complex}}, we have :
\begin{enumerate}
\item
For every $i=0,\dots,\dim X$ and $q\in\Z$, we have natural
isomorphisms of abelian groups :
$$
H_q(X^i,X^{i-1})\isom\left\{\begin{array}{ll}
                      \Z^{\Lambda_i} & \text{if $q=i$} \\
                      0            & \text{otherwise}.
                    \end{array}\right.
$$
\item
More precisely, let $(b_\lambda~|~\lambda\in\Lambda_i)$ be
the canonical basis of the free $\Z$-module $\Z^{\Lambda_i}$;
for every $\lambda\in\Lambda_i$, the isomorphism of\/ {\em(i)}
maps the image of the induced map
$$
H_if_\lambda:H_i(\bar\B{}^i(1),\SS^{i-1})\to H_i(X^i,X^{i-1})
$$
isomorphically onto the direct summand generated by $b_\lambda$.
\end{enumerate}
\end{lemma}
\begin{proof} (i): This is obvious for $i=0$. For
$i=1,\dots,k:=\dim X^\bullet$ and every $\lambda\in\Lambda_i$,
set
$$
Y^i_\lambda:=f_\lambda(\bar\B{}^i(1/2))
\qquad
Y^i:=\bigcup_{\lambda\in\Lambda_i}Y^i_\lambda
\qquad
A^i:=\bigcup_{\lambda\in\Lambda_i}\{f_\lambda(0)\}
$$
and for any $q\in\Z$, consider the natural group homomorphisms
\set\begin{equation}\label{eq_cell-complex}
{\diagram
H_q(\bar\B{}^i(1/2),\bar\B{}^i(1/2)\setminus\!\{0\})
\ar[r]^-{\alpha'} \ar[d]_\gamma &
H_q(\bar\B{}^i(1),\bar\B{}^i(1)\setminus\!\{0\}) \ar[d] &
\ar[l]_-{\beta'} \ar[d] H_q(\bar\B{}^i(1),\SS^{i-1}) \\
H_q(Y^i,Y^i\setminus A^i) \ar[r]^-\alpha &
H_q(X^i,X^i\setminus A^i) & \ar[l]_-\beta H_q(X^i,X^{i-1}).
\enddiagram}
\end{equation}
The map $\alpha$ is an isomorphism, by excision; the same
holds for $\beta$, since $X^{i-1}$ is a deformation retract
of $X^i\setminus A^i$. However, for every $\lambda\in\Lambda_i$
we have natural isomorphisms
$$
H_q(Y^i_\lambda,Y^i_\lambda\setminus A)\isom\left\{\begin{array}{ll}
                      \Z & \text{if $q=i$} \\
                      0  & \text{otherwise}.
                    \end{array}\right.
$$
whence the contention.

(ii): Take $q:=i$ in \eqref{eq_cell-complex}; then clearly $\gamma$
is a monomorphism whose image is the direct summand
$H_i(Y^i_\lambda,Y^i_\lambda\setminus A)$; on the other hand, arguing
as in the foregoing, we see that both $\alpha'$ and $\beta'$ are
isomorphisms. The assertion follows easily.
\end{proof}

\begin{remark}\label{rem_cell-complex}
(i)\ \ 
In the situation of lemma \ref{lem_cell-complex}, notice
that there is a natural bijection from $\Lambda_i$ to the set
of connected components $\pi_0(X^i\setminus X^{i-1})$ of
$X^i\setminus X^{i-1}$, for every $i=0,\dots,\dim X$.

(ii)\ \
The direct sum decomposition of $H_i(X^i,X^{i-1})$ provided by
lemma \ref{lem_cell-complex}(i) depends only on the filtration
$X^\bullet$ (and not on the choice of homeomorphisms $f_\lambda$).
Indeed, this is clear, since the map $H_if_\lambda$ factors
as a composition :
$$
H_i(\bar\B{}^i(1),\SS^{i-1})\isom
H_i(\bar e^i_\lambda,\bar e{}^i_\lambda\setminus e^i_\lambda)
\xrightarrow{\ l_\lambda\ }
H_i(X^i,X^{i-1})
$$
where $l_\lambda$ is the homomorphism induced by the obvious map
of pairs $(\bar e^i_\lambda,\bar e{}^i_\lambda\setminus e^i_\lambda)\to
(X^i,X^{i-1})$.

(iii)\ \
In view of (ii), the maps $l_\lambda$ admit natural left inverse
homomorphisms
$$
p_\lambda:H_i(X^i,X^{i-1})\to
H_i(\bar e^i_\lambda,\bar e{}^i_\lambda\setminus e^i_\lambda)
$$
such that
\set\begin{equation}\label{eq_decompose-homology}
p_\mu\circ l_\lambda=0
\qquad
\text{for every $\lambda,\mu\in\Lambda_i$ with $\lambda\neq\mu$}.
\end{equation}
These maps can be described as follows. For any
$\lambda\in\Lambda_i$, the natural map of pairs
$(X^i,X^{i-1})\to(X^i,X^i\setminus e^i_\lambda)$ induces a
homomorphism
$$
m_\lambda:H_i(X^i,X^{i-1})\to H_i(X^i,X^i\setminus e^i_\lambda).
$$
Arguing by excision as in the proof of lemma \ref{lem_cell-complex}(ii),
it is easily seen that $q_\lambda:=m_\lambda\circ l_\lambda$ is
an isomorphism (details left to the reader), and we set
$p_\lambda:=q^{-1}_\lambda\circ m_\lambda$. Obviously $p_\lambda$
is left inverse to $l_\lambda$, and in oder to check
\eqref{eq_decompose-homology} it suffices to show that
$m_\mu\circ l_\lambda=0$ for every $\mu\neq\lambda$.
However, $m_\mu\circ l_\lambda$ is the homomorphism arising
from the map of pairs
$(\bar e{}^i_\lambda,\bar e{}^i_\lambda\setminus e^i_\lambda)\to
(X^i,X^i\setminus e^i_\mu)$, so the assertion follows, after
remarking that $\bar e{}^i_\lambda\subset X^i\setminus e^i_\mu$.
\end{remark}

\sset\subsubsection{}\label{subsec_cell-complex}
We attach to any cell complex $X^\bullet$ a complex of abelian
groups $\cC_\bullet(X^\bullet)$, as follows.

$\bullet$\ \
For $i<0$ and for $i>k:=\dim X^\bullet$, we set $\cC_i(X^\bullet):=0$,
and for $i=1,\dots,k$ we let
$$
\cC_i(X^\bullet):=H_i(X^i,X^{i-1}).
$$
Sometimes we write just $\cC_i(X)$, unless it is useful to
stress which filtration $X^\bullet$ on $X$ we are considering. 
For every $i>0$, the differential
$d_i:\cC_i(X^\bullet)\to\cC_{i-1}(X^\bullet)$ is the composition
$$
H_i(X^i,X^{i-1})\xrightarrow{\ \partial_i\ }
H_{i-1}(X^{i-1})\xrightarrow{\ j_{i-1}\ } H_{i-1}(X^{i-1},X^{i-2})
$$
where $\partial_i$ is the boundary operator of the long exact
homology sequence associated with the pair $(X^i,X^{i-1})$, and
$j_{i-1}$ is the homomorphism induced by the obvious map of pairs
$(X^{i-1},\emptyset)\to(X^{i-1},X^{i-2})$. In order to check that
$d_i\circ d_{i+1}=0$ for every $i\in\Z$, recall that every
element of $H_{i+1}(X^{i+1},X^i)$ is the class $\bar c$ of a singular
$(i+1)$-chain $c$ of $X^{i+1}$, whose boundary $\partial_{i+1}c$
is a singular $i$-chain of $X^i$; then $d_{i+1}(\bar c)$ is the class
of $\partial_{i+1}c$ in $H_i(X^i,X^{i-1})$, and $d_i\circ d_{i+1}(\bar c)$
is the class of $\partial_i\circ\partial_{i+1}c$ in $H_{i-1}(X^{i-1},X^{i-2})$,
so it vanishes.

$\bullet$\ \
It is also useful to consider an augmented version of the above
complex; namely, let us set
$$
\bar\cC_{-1}(X^\bullet):=\Z
\qquad\text{and}\qquad
\bar\cC_i(X^\bullet):=\cC_i(X^\bullet)
\qquad
\text{for every $i\neq -1$}
$$
with differential $d_0$ given by the rule : $d_0(b^0_\lambda)=1$
for every $\lambda\in\Lambda_0$.

\begin{proposition}\label{prop_computes-sing-hom}
With the notation of \eqref{subsec_cell-complex},
there exist natural isomorphisms of abelian groups
$$
H_q\cC_\bullet(X^\bullet)\isom H_q(X)
\qquad
\text{for every $q\in\N$}.
$$
\end{proposition}
\begin{proof} For every topological space $T$, let $C_\bullet(T)$
denote the complex of singular chains of $T$. The filtration
$X^\bullet$ induces a finite filtration of complexes
$$
C_\bullet(X^0)\subset C_\bullet(X^1)\subset C_\bullet(X^2)
\subset\cdots\subset C_\bullet(X)
$$
whence a convergent spectral sequence
$$
E^1_{pq}:=H_{p+q}(X^p,X^{p-1})\Rightarrow H_{p+q}(X)
$$
(see \cite[Th.5.5.1]{We}). By direct inspection, it is easily
seen that the differential $d^1_{p,0}:E^1_{p,0}\to E^1_{p-1,0}$
agrees with the differential $d_p$ of $\cC_\bullet(X^\bullet)$,
for every $p\in\N$. On the other hand, lemma \ref{lem_cell-complex}
shows that this spectral sequence degenerates, whence the
contention.
\end{proof}

\sset\subsubsection{}\label{subsec_incidence}
Keep the notation of definition \ref{def_cell-complex}, and
set $\Lambda_0:=X^0$. For every $i=0,\dots,\dim X$ and every
$\lambda\in\Lambda_i$, the abelian group
$H_{i,\lambda}:=H_i(\bar e^i_\lambda,\bar e{}^i_\lambda\setminus e^i_\lambda)$
is free of rank one; if $i>0$, we fix one of the two generators
of this group, and we denote by $b^i_\lambda$ its image in
$\cC_i(X^\bullet)$ (see remark \ref{rem_cell-complex}(ii)).
For $i=0$, each $e^0_\lambda$ is a point, hence $H_{0,\lambda}$
admits a canonical identification with $\Z$, and we let $b^0_\lambda$
be the image of $1$, under the resulting map
$\Z\isom H_{0,\lambda}\to\cC_0(X^\bullet)$. The system of classes
$(b^i_\lambda~|~0\leq i\leq\dim X,\lambda\in\Lambda_i)$ is called
an {\em orientation\/} for $X^\bullet$. We may write
$$
d_i(b^i_\lambda)=\sum_{\mu\in\Lambda_{i-1}}[b^i_\lambda:b^{i-1}_\mu]b^{i-1}_\mu
\qquad
\text{for every $i\leq\dim X^\bullet$ and every $\lambda\in\Lambda_i$}
$$
for a system of uniquely determined integers $[b^i_\lambda:b^{i-1}_\mu]$,
called the {\em incidence numbers} of the cells $e^i_\lambda$ and
$e^{i-1}_\mu$ (relative to the chosen orientation of $X^\bullet$).
With this notation, we may state :

\begin{lemma}\label{lem_incidence}
The incidence numbers of a cell complex $X^\bullet$ fulfill the
following conditions :
\begin{enumerate}
\item
$\sum_{\mu\in\Lambda_{i-1}}
[b^i_\lambda:b^{i-1}_\mu]\cdot[b^{i-1}_\mu:b^{i-2}_\nu]=0$ for every
$i\geq 2$, every $\lambda\in\Lambda_i$ and every $\nu\in\Lambda_{i-2}$.
\item
Let $\lambda\in\Lambda_1$ be any index, and say that
$\bar e{}^1_\lambda\setminus e^1_\lambda=e^0_\mu\cup e^0_\rho$.
Then $[b^1_\lambda:b^0_\mu]+[b^1_\lambda:b^0_\rho]=0$.
\end{enumerate}
\end{lemma}
\begin{proof} Condition (i) translates the identity
$d_{i-1}\circ d_i=0$ for the differential $d_\bullet$ of
the complex $\cC_\bullet(X^\bullet)$. The identity of (ii)
follows by a simple inspection of the definition of $d_1$ :
details left to the reader.
\end{proof}

\sset\subsubsection{}\label{subsec_cone-cell}
The generalities of the previous paragraphs shall be applied
to the following situation. Let $(V,\sigma)$ be a strictly
convex polyhedral cone such that $\La\sigma\Ra=V$, and set
$d:=\dim_\R V$. We attach to $\sigma$ a cell complex
$C_\sigma^\bullet$ as follows. Fix $u_0\in\sigma^\vee$ such that
$\sigma\cap\Ker\,u_0=\{0\}$ (corollary \ref{cor_strongly}),
let $\sigma^\circ$ be the topological interior of $\sigma$
(in $V$) and set
$$
C_\sigma:=\sigma\cap u_0^{-1}(1)
\qquad
C^\circ_\sigma:=\sigma^\circ\cap u^{-1}_0(1).
$$
Let $v_1,\dots,v_k$ be a minimal system of generators for $\sigma$;
we may assume that $u_0(v_i)=1$ for $i=1,\dots,k$, in which case
\set\begin{equation}\label{eq_image-of-simplex}
C_\sigma=\biggl\{\sum_{i=1}^k\lambda_iv_i~|~\lambda_1,\dots,\lambda_k\geq 0
\quad\text{and}\quad\sum_{i=1}^k\lambda_i=1\biggr\}
\end{equation}
which shows that $C_\sigma$ is a compact topological space
(with the topology inherited from $V$).

\begin{lemma}\label{lem_cone-cell}
There exists a homeomorphism $\bar\B{}^{d-1}(1)\isom C_\sigma$
that maps $\B^{d-1}(1)$ onto $C^\circ_\sigma$.
\end{lemma}
\begin{proof} Set $v_0:=k^{-1}\cdot(v_1+\cdots+v_k)$ and
$W:=\Ker\,u_0$. It suffices to show that there exists a
homeomorphism $\bar\B{}^{d-1}(1)\isom D_\sigma:=C_\sigma-v_0\subset W$
that maps $\B{}^{d-1}(1)$ onto $D^\circ_\sigma:=C_\sigma^\circ-v_0$.
However, pick a minimal system $u_1,\dots,u_t$ of generators of
$\sigma^\vee$ (corollary \ref{cor_spanning-cone}(i)); since
$\sigma$ spans $V$, for every $i=1,\dots,t$ there exists $j\leq k$
such that $u_i(v_j)>0$. Hence -- after replacing the $u_i$ by
suitable scalar multiples -- we may assume that $u_i(v_0)=1$ for
every $i=1,\dots,t$, and then lemma \ref{lem_ddoble-cone} implies that
$$
D_\sigma=
\{v\in W~|~u_i(v)\geq-1\quad\text{for every $i=1,\dots,t$}\}.
$$
Also, proposition \ref{prop_spanning-cone}(i) implies that
$$
D^\circ_\sigma=\{v\in W~|~u_i(v)>-1\quad\text{for every $i=1,\dots,t$}\}.
$$
For every $v\in W$, set $\mu(v):=\min(u_i(v)~|i=1,\dots,t)$. It
is easily seen that $\mu(v)<0$ for every $v\in W\setminus\{0\}$;
indeed, if $\mu(w)\geq 0$, then
the subset $\{\lambda w~|~\lambda\in\R_+\}$ lies in $D_\sigma$, and
since the latter is compact, it follows that $w=0$. Moreover we have
\set\begin{equation}\label{eq_segments}
(\R_+w)\cap D_\sigma=\{\lambda w~|~\lambda\in[0,-1/\mu(w)]\}
\qquad
\text{for every $w\in W\setminus\{0\}$}.
\end{equation}
Fix a Banach norm $\Vert\cdot\Vert_V$ on $V$, and consider the
mapping
$$
\phi:D_\sigma\to W
\qquad\text{such that}\qquad
\phi(0):=0
\quad\text{and}\quad
\phi(v):=-\frac{\mu(v)}{\Vert v\Vert}\cdot v
\quad
\text{for $v\neq 0$}.
$$
It is easily seen that $\phi$ is injective; since $\mu$ is
a continuous mapping, the same follows for $\phi$, and since
$D_\sigma$ is compact, we conclude that $\phi$ induces a
homeomorphism $D_\sigma\isom\phi(D_\sigma)$. However, from
\eqref{eq_segments} it follows that
$$
\phi(D_\sigma)=\{v\in W~|~\Vert v\Vert\leq 1\}
\qquad\text{and}\qquad
\phi(D^\circ_\sigma)=\{v\in W~|~\Vert v\Vert<1\}
$$
whence the claim.
\end{proof}

\sset\subsubsection{}\label{subsec_filter-cone}
Keep the notation of \eqref{subsec_cone-cell}; we consider the
finite filtration $C^\bullet_\sigma$ of $C_\sigma$, defined as
follows. For every $i=-1,\dots,d-1$, we let $C^i_\sigma\subset C_\sigma$
be the union of the subsets $\tau\cap u_0^{-1}(1)$, where $\tau$
ranges over the (finite) set of all faces of $\sigma$ of dimension
$i+1$. Clearly $C^{-1}_\sigma=\emptyset$, $C^{d-1}_\sigma=C_\sigma$,
and $C^i_\sigma$ is a closed subset of $C_\sigma$, for every
$i=0,\dots,d-1$. Moreover, it follows easily from lemma
\ref{lem_cone-cell} and proposition \ref{prop_spanning-cone}(i)
that the datum of $C_\sigma$ and its filtration $C^\bullet_\sigma$
is a cell complex. With the notation of definition
\ref{def_cell-complex}, the indexing set $\Lambda_i$ can be
taken to be the set of all $(i+1)$-dimensional faces of $\sigma$,
for every $i=0,\dots,d-1$ : indeed, if $\tau$ is such a face,
denote by $\tau^\circ$ the relative interior of $\tau$ (see
example \ref{ex_cone-dim-two}(iii)); then it is clear that
$e^i_\tau:=\tau^\circ\cap C_\sigma\neq\emptyset$ and
$\bar e{}^i_\tau=\tau\cap C_\sigma$.

\begin{remark}\label{rem_cone-complex}
(i)\ \
Notice that the cell complex $C^\bullet_\sigma$ is independent
-- up to homeomorphism -- of the choice of $u_0$. Indeed, say
that $u'_0\in\sigma^\vee$ is any other linear form such that
$\sigma\cap\Ker\,u'_0=\{0\}$, and let
$C'_\sigma:=\sigma\cap u^{\prime-1}_0(1)$. We define a homeomorphism
$\psi:C'_\sigma\isom C_\sigma$, by the rule :
$$
v\mapsto u_0(v)^{-1}\cdot v
\qquad
\text{for every $v\in C'_\sigma$}.
$$
It is easily seen that $\psi$ restricts to homeomorphisms
$C^{\prime i}_\sigma\isom C^i_\sigma$ for every $i=0,\dots,d-1$.

(ii)\ \
For any integer $i=0,\dots,d-3$, and any cells $e^i_\tau$,
$e^{i+2}_\lambda$ of $C^\bullet_\sigma$ such that
$e^i_\tau\subset\bar e{}^{i+2}_\lambda$, there exist exactly
two $(i+1)$-dimensional cells $e^{i+1}_\mu$, $e^{i+1}_\rho$ such
that $e^i_\tau\subset\bar e{}^{i+2}_\mu\cap\bar e^{i+1}_\rho$ and
$e^{i+1}_\mu\cup e^{i+1}_\rho\subset\bar e{}^{i+2}_\lambda$ :
indeed, this assertion is a direct translation of claim
\ref{cl_codim-face}(ii).
\end{remark}

\begin{proposition}\label{prop_incidence-nos}
With the notation of \eqref{subsec_incidence} and
\eqref{subsec_filter-cone}, the following holds for every
$i=0,\dots,d-2$ :
\begin{enumerate}
\item
If\/ $\tau\in\Lambda_{i+1}$ and $\mu\in\Lambda_i$ is not
a facet of\/ $\tau$, we have $[b_\tau^{i+1}:b^i_\mu]=0$.
\item
If\/ $\tau\in\Lambda_{i+1}$ and $\mu\in\Lambda_i$ is a facet of\/
$\tau$, then $[b^{i+1}_\tau:b^i_\mu]\in\{1,-1\}$.
\item
In the situation of remark {\em\ref{rem_cone-complex}(ii)}, we have :
$$
[b^{i+2}_\lambda:b^{i+1}_\mu]\cdot[b^{i+1}_\mu:b^i_\tau]+
[b^{i+2}_\lambda:b^{i+1}_\rho]\cdot[b^{i+1}_\rho:b^i_\tau]=0.
$$
\end{enumerate}
\end{proposition}
\begin{proof}(i): Consider the commutative diagram
\set\begin{equation}\label{eq_calculate-incidence}
{\diagram
H_{i+1}(\bar e{}^{i+1}_\tau,\bar e{}^{i+1}_\tau\setminus e^{i+1}_\tau)
\ar[r]^-\partial \ar[d]_{l_\tau} &
H_i(\bar e{}^{i+1}_\tau\setminus e^{i+1}_\tau) \ar[r]^-{j_\tau} \ar[d]^-g &
H_i(C^i_\sigma,C^i_\sigma\setminus e^i_\mu) \\
H_{i+1}(C^{i+1}_\sigma,C^i_\sigma) \ar[r]^-{\partial'} & H_i(C^i_\sigma)
\ar[r]^-{j_i} & H_i(C^i_\sigma,C^{i-1}_\sigma) \ar[u]_{m_\mu}
\enddiagram}
\end{equation}
where $l_\tau$ and $m_\mu$ are defined as in remark
\ref{rem_cell-complex}(ii,iii) and $g$ is induced by the inclusion
map $\bar e{}^{i+1}_\tau\setminus e^{i+1}_\tau\to C^i_\sigma$, the
maps $\partial$ and $\partial'$ are the boundary operators of the
long exact sequences attached to the pairs
$(\bar e{}^{i+1}_\tau,\bar e{}^{i+1}_\tau\setminus e^{i+1}_\tau)$ and
$(C^{i+1}_\sigma,C^i_\sigma)$, and $j_i$ (resp. $j_\tau$) is
deduced from the obvious map of pairs
$(\bar e{}^{i+1}_\tau\setminus e^{i+1}_\tau,\emptyset)\to
(C^i_\sigma,C^i_\sigma\setminus e^i_\mu)$
(resp. $(C^i_\sigma,\emptyset)\to(C^i_\sigma,C^{i-1}_\sigma)$).
In light of remark \ref{rem_cell-complex}(iii), we need to
check that $m_\mu\circ j_{i-1}\circ\partial'\circ l_\tau=0$,
and to this aim, it suffices to show that $j_\tau=0$.
However, the assumption implies that $\mu\cap\tau$ is a (proper)
face of both $\mu$ and $\tau$; this translates as the identity
$(\bar e{}^{i+1}_\tau\setminus e^{i+1}_\tau)\cap e^i_\mu=\emptyset$,
whence the claim.

(ii): Clearly $C^\bullet_\tau$ is a cell subcomplex of $C^\bullet_\sigma$,
and $e^{i+1}_\tau$ and $e^i_\mu$ are cells of this subcomplex.
Moreover, any orientation for $C^\bullet_\sigma$ restricts to
an orientation for $C^\bullet_\tau$, and the resulting incidence
number of $e^{i+1}_\tau$ and $e^i_\mu$ is the same for either cell
complex. Thus, we may replace $\sigma$ by $\tau$, in which case
the map $g$ of \eqref{eq_calculate-incidence} is the identity.

\begin{claim}\label{cl_should-reduce}
If $\sigma=\tau$ and $i>0$, both $\partial$ and $k_\mu:=m_\mu\circ j_i$
in \eqref{eq_calculate-incidence} are isomorphisms.
\end{claim}
\begin{pfclaim} For $\partial$, we use the long exact homology
sequence of the pair
$(\bar e{}^{i+1}_\tau,\bar e{}^{i+1}_\tau\setminus e^{i+1}_\tau)$:
since $i\geq 0$, lemma \ref{lem_cone-cell} implies that
$H_{i+1}(\bar e^{i+1}_\tau)=0$, so $\partial$ is injective;
since $i>0$, the same argument shows that $\partial$ is surjective.
Next, $k_\mu$ is the homomorphism deduced from the long
exact homology sequence of the pair
$(C^i_\tau,C^i_\tau\setminus e^i_\mu)$, and we remark that
$C^i_\tau\setminus e^i_\mu$ is contractible : indeed, $e^i_\mu$
is homeomorphic to $\B^i(1)$ (lemma \ref{lem_cone-cell}), hence
$C^i_\tau\setminus e^i_\mu$ is a retraction of $C^i_\tau\setminus\{x\}$,
for any $x\in e^i_\mu$, and the latter is homeomorphic to $\R^i$
(again by lemma \ref{lem_cone-cell}). If $i>1$, we deduce already
that $k_\mu$ is an isomorphism. For $i=1$, the same argument shows
that $k_\mu$ is injective, so its cokernel is a cyclic torsion
group that injects into $H_0(C^1_\tau\setminus e^1_\mu)\simeq\Z$,
hence it must vanish as well.
\end{pfclaim}

Since \eqref{eq_calculate-incidence} commutes, claim
\ref{cl_should-reduce} yields the assertion, in case $i>0$.
If $i=0$, $C^1_\tau$ is isomorphic to $\bar\B{}^1(1)$, and
$e^0_\mu$ is one of the two points of
$\bar\B{}^1(1)\setminus\B^1(1)$, so the assertion can be checked
easily, by inspecting the definitions.

(iii) follows from (i), lemma \ref{lem_incidence}(i) and remark
\ref{rem_cone-complex}(ii).
\end{proof}

The following result says that the properties of proposition
\ref{prop_incidence-nos} completely characterize the incidence
numbers of $C^\bullet_\sigma$.

\begin{proposition} Keep the notation of \eqref{subsec_filter-cone},
and consider a system of integers
$$
(\beta^i_{\lambda\mu}~|~
i=1,\dots,d-1,\ \lambda\in\Lambda_i,\mu\in\Lambda_{i-1})
$$
such that :
\begin{enumerate}
\alphaenu
\item
If\/ $\lambda\in\Lambda_i$ and $\mu\in\Lambda_{i-1}$ is not
a facet of\/ $\lambda$, then $\beta^i_{\lambda\mu}=0$.
\item
If\/ $\lambda\in\Lambda_i$ and $\mu\in\Lambda_{i-1}$ is a
facet of\/ $\lambda$, then $\beta^i_{\lambda\mu}\in\{1,-1\}$.
\item
If $\mu$ and $\tau$ are the two $1$-dimensional facets of
the $2$-dimensional face $\lambda$ of $\sigma$, then
$$
\beta^1_{\lambda\mu}+\beta^1_{\lambda\tau}=0.
$$
\item
Let $\lambda\in\Lambda_{i+1}$ and $\mu\in\Lambda_{i-1}$ be two faces,
such that $\mu$ is a face of $\lambda$, and denote by $\tau$ and
$\rho$ the two facets of $\lambda$ that contain $\mu$. Then
$$
\beta^{i+1}_{\lambda\tau}\beta^i_{\tau\mu}+
\beta^{i+1}_{\lambda\rho}\beta^i_{\rho\mu}=0.
$$
\end{enumerate}
Then there exists a unique orientation of $C^\bullet_\sigma$
$$
(b^i_\lambda~|~i=0,\dots,d-1,\ \lambda\in\Lambda_i)
$$
such that for every $i=1,\dots,d-1$ we have :
\set\begin{equation}\label{eq_orientation}
[b^i_\lambda:b^{i-1}_\mu]=\beta^i_{\lambda\mu}
\qquad
\text{for every $\lambda\in\Lambda_i$, $\mu\in\Lambda_{i-1}$}.
\end{equation}
\end{proposition}
\begin{proof} We construct the orientation classes $b^i_\lambda$
fulfilling condition \eqref{eq_orientation}, by induction on $i$.
For $i=0$, the condition is empty, and the classes $b^0_\lambda$
are prescribed by \eqref{subsec_incidence}. For $i=1$, and a
given $\lambda\in\Lambda_1$, denote by $\tau$ and $\rho$ the
two faces of $\lambda$; clearly the condition
$[b^1_\lambda:b^0_\tau]=\beta^1_{\lambda\tau}$ (and assumption (b))
determines $b^1_\lambda$ univocally, and in view of (c) and lemma
\ref{lem_incidence}(ii), for this choice of orientation of
$e^1_\lambda$ we have as well
$[b^1_\lambda:b^0_\rho]=\beta^1_{\lambda\rho}$. Moreover, if
$\mu\in\Lambda_0\setminus\{\tau,\rho\}$, then \eqref{eq_orientation}
is verified as well, by virtue of (a) and proposition
\ref{prop_incidence-nos}(i).

Now, suppose that $i>1$ and that we have already constructed
orientation classes $b^j_\mu$ as sought, for every $j<i$ and
every $\mu\in\Lambda_j$. Let $\lambda\in\Lambda_i$ be any face,
and fix a facet $\tau$ of $\lambda$; again, the identity
\set\begin{equation}\label{eq_fix-lambda-tau}
[b^i_\lambda:b^i_\tau]=\beta^i_{\lambda\tau}
\end{equation}
determines $b^i_\lambda$, and it remains to check that -- with
this choice of $b^i_\lambda$ -- condition \eqref{eq_orientation}
holds for every $\mu\in\Lambda_{i-1}\setminus\{\tau\}$, and
by virtue of (a) and proposition \ref{prop_incidence-nos}(i),
it suffices to consider the facets $\mu$ of $\lambda$.
However, notice that the system of orientation classes
$(b^j_\mu~|~j=0,\dots,i-1,\ \mu\subset\lambda)$ amounts to
an orientation of $C^{i-1}_\lambda$ (regarded as a cell
subcomplex of $C^\bullet_\sigma$), and the incidence numbers
for the complex $(C^{i-1}_\lambda)^\bullet$ relative to these
classes agree with the incidence numbers of $C^\bullet_\sigma$,
relative to the same classes. Especially, the sums
$$
c:=\sum_{\mu\in\Lambda_{i-1}}\beta^i_{\lambda\mu}b^{i-1}_\mu
\qquad
c':=\sum_{\mu\in\Lambda_{i-1}}[b^i_\lambda:b^{i-1}_\mu]b^{i-1}_\mu
$$
are well defined elements of $\cC_{i-1}(C^{i-1}_\lambda)$, and
in fact :
$$
\begin{aligned}
d_{i-1}(c) &
=\sum_{\mu\in\Lambda_{i-1}}\beta^i_{\lambda\mu}\cdot d_{i-1}(b^{i-1}_\mu) \\
& =\sum_{\mu\in\Lambda_{i-1}}\beta^i_{\lambda\mu}\cdot\sum_{\rho\in\Lambda_{i-2}}
[b^{i-1}_\mu:b^{i-2}_\rho]b^{i-2}_\rho \\
& =\sum_{\mu\in\Lambda_{i-1}}\sum_{\rho\in\Lambda_{i-2}}
\beta^i_{\lambda\mu}\beta^{i-1}_{\mu\rho}b^{i-2}_\rho &
\text{(by inductive assumption)} \\
& =0 & \text{(by (d))}
\end{aligned}
$$
and a similar calculation yields $d_{i-1}(c')=0$ as well.
However, $C^{i-1}_\lambda$ is homeomorphic to $\SS^{i-1}$
(lemma \ref{lem_cone-cell}), hence
$$
\Ker(d_{i-1}:\cC_{i-1}(C^{i-1}_\lambda)\to \cC_{i-2}(C^{i-1}_\lambda))
\simeq\Z
$$
which, in view of (b), implies that $c=\pm c'$. Taking into
account \eqref{eq_fix-lambda-tau}, we see that actually
$c=c'$, whence the claim. The uniqueness of the orientation
fulfilling condition \eqref{eq_orientation} can be checked
easily by induction on $i$ : the details shall be left to
the reader.
\end{proof}

\sset\subsubsection{}\label{subsec_back-to-cohomology}
Let now $L$ be a free abelian group of finite rank $d$, and
$(L_\R,\sigma)$ a strictly convex $L$-rational polyhedral
cone (see \eqref{subsec_from-con-to-mon}), such that
$\La\sigma\Ra=L_\R$; set
$$
P:=L\cap\sigma
\qquad
F_\lambda:=L\cap\lambda
\qquad
P_\lambda:=F_\lambda^{-1}P
\qquad
\text{for every face $\lambda$ of $\sigma$}
$$
so $P$ and its localization $P_\lambda$ are fine and saturated
monoids (proposition \ref{prop_Gordon}(i)), and $F_\lambda$
is a face of $P$, for every such $\lambda$. Let $R$ be any
ring, and if $\lambda,\mu\subset\sigma$ are any two faces with
$\lambda\subset\mu$, denote by
$$
j_{\lambda\mu}:R[P_\lambda]\to R[P_\mu]
$$
the natural localization map; notice that if $0\subset\sigma$
is the unique $0$-dimensional face, then $P_0=P$, and $j_{0\mu}$
is the localization map $R[P]\to R[P_\mu]$. We attach to $P$
the complex $\bar\cC{}^\bullet_P$ of $R[P]$-modules such that :
$$
\bar\cC{}^0_P:=R[P]
\qquad\text{and}\qquad
\bar\cC{}^i_P:=\bigoplus_{\lambda\in\Lambda_{i-1}}R[P_\lambda]
\qquad
\text{for every $i=1,\dots,\dim P$}
$$
with differentials given by the rule :
$$
d^i(x_\lambda):=
\sum_{\substack{\mu\in\Lambda_i \\
               \lambda\subset\mu}}
[b^i_\mu:b^{i-1}_\lambda]\cdot j_{\lambda\mu}(x_\lambda)
\qquad
\text{for every $i\geq 1$, $\lambda\in\Lambda_{i-1}$ and
$x_\lambda\in R[P_\lambda]$}
$$
and
$$
d^0(x):=\sum_{\mu\in\Lambda_0}j_{0\mu}(x)
\qquad
\text{for every $x\in R[P]$}
$$
where $(b^j_\tau~|~j=0,\dots\dim\sigma-1,\ \tau\in\Lambda_j)$
is a chosen orientation for $C^\bullet_\sigma$. Taking into
account lemma \ref{lem_incidence}(ii) and proposition
\ref{prop_incidence-nos}(iii), it is easily seen that
$d^{i+1}\circ d^i=0$ for every $i\geq 0$ (essentially,
the differential $d^i$ is the transpose of the differential
$d_{i-1}$ of $\bar\cC_\bullet(C^\bullet_\sigma)$).

Set as well $X:=\Spec\,R[P]$, and notice that
$Z:=\Spec\,R\La P/\fm_P\Ra\simeq\Spec\,R$ is a closed subset
of $X$. For any $R[P]$-module $M$, we denote as usual by
$M^\sim$ the quasi-coherent $\cO_{\!X}$-module arising from $M$.

\begin{theorem}\label{th_combinatorial-coh-compute}
With the notation of \eqref{subsec_back-to-cohomology}, for
every $R[P]$-module $M$ we have a natural isomorphism
$$
R\Gamma_{\!Z}M^\sim\isom M\otimes_{R[P]}\bar\cC{}^\bullet_P
\qquad
\text{in $\sD^+(R[P]\Mod)$}.
$$
\end{theorem}
\begin{proof} For every face $\lambda$ of $\sigma$, set
$U_\lambda:=\Spec\,R[P_\lambda]$, and denote by
$g_\lambda:U_\lambda\to X$ the open immersion. We consider
the chain complex $\cR_\bullet$ of $\Z_X$-modules such that :
$$
\cR_0:=\Z_X
\qquad
\cR_i:=\bigoplus_{\lambda\in\Lambda_{i-1}}g_{\lambda!}\Z_{U_\lambda}
\qquad
\text{for every $i=1,\dots,\dim\,P$}.
$$
The differential $d_1:\cR_1\to\cR_0$ is just the sum of the
natural morphisms $g_{\lambda!}\Z_{U_\lambda}\to\Z_X$, for $\lambda$
ranging over the one-dimensional faces of $\sigma$. For $i>1$,
the differential $d_i$ is the sum of the maps
$$
d_\lambda:=\sum_{\substack{\mu\in\Lambda_i \\
               \mu\subset\lambda}}
[b^i_\lambda:b^{i-1}_\mu]\cdot d_{\lambda\mu!}:
g_{\lambda!}\Z_{U_\lambda}\to\cR_{n-1}
\qquad
\text{for every $\lambda\in\Lambda_{i-1}$}
$$
where $d_{\lambda\mu}:
g_{\lambda!}\Z_{U_\lambda}\to g_{\mu!}\Z_{U_\mu}\subset\cR_{n-1}$
is induced by the inclusion $U_\lambda\subset U_\mu$. With
this notation, a simple inspection of the definitions yields
a natural identification
\set\begin{equation}\label{eq_inspect-this}
M\otimes_{R[P]}\bar\cC{}^\bullet_P\isom
\Hom^\bullet_\Z(\cR_\bullet,M^\sim[0])
\qquad
\text{in $\sC(R[P]\Mod)$}.
\end{equation}
On the other hand, let $\lambda_1,\dots,\lambda_n$ be the
one-dimensional faces of $\sigma$; then $\lambda_i$ is
$L$-rational (see \eqref{subsec_from-con-to-mon}), and
$F_{\lambda_i}$ is a fine, sharp and saturated monoid of
dimension one (propositions \ref{prop_Gordon}(i) and
\ref{prop_Gordie}(ii)), so $F_{\lambda_i}\simeq\N$ for $i=1,\dots,n$
(theorem \ref{th_structure-of-satu}(ii)). For each $i=1,\dots,n$,
let $y_i$ be the unique generator of $F_{\lambda_i}$, and denote
by $I\subset P$ the ideal generated by $y_1,\dots,y_n$; we remark

\begin{claim}\label{cl_radical-of-I}
The radical of $I$ is $\fm_P$.
\end{claim}
\begin{pfclaim} For any $x\in\fm_P$, pick a subset
$S\subset\{1,\dots,n\}$ such that $x=\sum_{i\in S}a_iy_i$,
with $a_i>0$ for every $i\in S$. Let $N\in\N$ be large enough
such that $Na_i\geq 1$ for every $i\in S$, and pick $b_i\in\N$
such that $Na_i\geq b_i\geq 1$ for every $i\in S$. It follows
that $Nx-\sum_{i\in S}b_iy_i\in P$, and the claim follows.
\end{pfclaim}

Now, let $i:Z\to X$ be the closed immersion; we have :

\begin{claim}\label{cl_combinatorial-Cech}
The natural map $\Z_X\to i_*\Z_Z$ induces an isomorphism
\set\begin{equation}\label{eq_new-cech-resolve}
\cR_\bullet\isom i_*\Z_Z[0]
\qquad
\text{in $\sD^+(\Z_X\Mod)$}.
\end{equation}
\end{claim}
\begin{pfclaim} The assertion can be checked on the stalks
at the points of $X$, hence let $x\in X$ be any such point;
if $x\in Z$, claim \ref{cl_radical-of-I} easily implies
that $(\cR_\bullet)_x$ is concentrated in degree zero, and
then $\eqref{eq_new-cech-resolve}_x$ is clearly an isomorphism
of complexes of $\Z$-modules. If $x\notin Z$, let $\fp\subset R[P]$
be the prime ideal corresponding to $x$, and
$\lambda\subset\sigma$ the unique face such that
$P\setminus F_\lambda=\fp\cap P$; a simple inspection
of the construction shows that
$$
(\cR_\bullet)_x\isom\bar\cC_\bullet(C_\lambda^\bullet)
$$
(notation of \eqref{subsec_cell-complex}). But proposition
\ref{prop_computes-sing-hom} and lemma \ref{lem_cone-cell}
imply that $\bar\cC_\bullet(C^\bullet_\lambda)$ is acyclic,
whence the claim.
\end{pfclaim}

In view of \eqref{eq_inspect-this} and claim
\ref{cl_combinatorial-Cech}, we are reduced to showing that
the natural map
$$
\Hom^\bullet_\Z(\cR_\bullet,M^\sim[0])\to
R\Hom^\bullet_\Z(\cR_\bullet,M^\sim[0])
$$
is an isomorphism in $\sD^+(R[P]\Mod)$. This can be done by
a spectral sequence argument, along the lines of proposition
\ref{prop_affine-Cech-resolve} (which indeed includes a special
case of the situation we are considering here, namely the
case where $P$ is a free monoid : the details shall be left
to the reader).
\end{proof}

\sset\subsubsection{}
Notice that $\bar\cC{}^\bullet_P$ is a complex of $L$-graded
$R$-modules, hence its cohomology is $L$-graded as well,
and we wish next to compute the graded terms
$\gr_\bullet H^\bullet(\bar\cC{}^\bullet_P)$. To this aim, we
make the following :

\begin{definition} Let $V$ be a finite dimensional $\R$-vector
space, $X\subset V$ any subset, and $z\in V$ any point. Then :
\begin{enumerate}
\alphaenu
\item
We say that a point $x\in X$ {\em is visibile from $z$},
if $\{tz+(1-t)x~|~0\leq t\leq 1\}\cap X=\{x\}$.
\item
We say that a subset $S\subset X$ {\em is visibile from $z$},
if every point of $S$ is visible from $z$.
\end{enumerate}
\end{definition}

\begin{lemma}\label{lem_do-you-see-me}
In the situation of \eqref{subsec_filter-cone}, let
$z\in V\setminus C_\sigma$ be any point, and denote by
$S$ the set of points of\/ $C_\sigma$ that are visible from $z$.
Then :
\begin{enumerate}
\item
$S$ is a subcomplex of\/ $C^\bullet_\sigma$.
\item
$S$ (with the topology induced from $V$) is homeomorphic to
$\bar\B{}^e(1)$, for $e\in\{d-2,d-1\}$.
\end{enumerate}
\end{lemma}
\begin{proof} Pick a system of generators $\rho_1,\dots,\rho_n$
for $\sigma^\vee$, and choose $\rho_0\in\sigma^\vee$ so that
$C_\sigma=\sigma\cap\rho_0^{-1}(1)$; we may assume that, for
some integer $k\leq n$ we have
$$
\rho_i(z)<0
\qquad\text{if and only if}\qquad
1\leq i\leq k.
$$
Say that $y\in C_\sigma$, so that $\rho_i(y)\geq 0$ for every
$i=1,\dots,n$, and set $y_t:=tz+(1-t)y$ for every $t\in[0,1]$;
then clearly $\rho_i(y_t)\geq 0$ for every $i=k+1,\dots,n$.
Now, if $\rho_0(z)\neq 1$, we get $\rho_0(y_t)\neq 1$ for every
$t\neq 0$, therefore the whole of $C_\sigma$ is visible from $z$,
in which case (i) is trivial, and (ii) follows from \ref{lem_cone-cell}.
Hence we may assume that $\rho_0(z)=1$, so $\rho_0(y_t)=1$ for
every $t\in[0,1]$.
Suppose now that $\rho_i(y)>0$ for some $i\leq k$; then there
exists a unique $t_i\in]0,1[$ such that $\rho_i(y_{t_i})=0$.
Hence, if $s:=\min(\rho_i(y)~|~i=1,\dots,k)>0$, let
$t:=\min(t_i~|~i=1,\dots,k)$; it follows that $\rho_i(y_t)\geq 0$
for every $i=1,\dots,n$, so $y_t\in C_\sigma$, which says that $y$
is not visible from $z$. Conversely, if $s=0$, then it follows
easily that $y$ is visible from $z$. We conclude that
$$
S=C_\sigma\cap\bigcup_{i=1}^k\Ker\,\rho_i
$$
which shows (i). Next, set $W:=\Ker\,\rho_0$, and denote by
$\tau_z:\rho_0^{-1}(1)\isom W$ the translation map given by
the rule : $x\mapsto x-z$ for every $x\in\rho_0^{-1}(1)$.
To conclude, it suffices to check that $S':=\tau_z(S)$ is
homeomorphic to $\bar\B{}^{d-2}(1)$. To this aim, denote by
$\lambda$ the convex cone in $W$ generated by $\tau_z(C_\sigma)$.
Explicitly, if $v_1,\dots,v_k\in\sigma$ have been chosen so that
\eqref{eq_image-of-simplex} holds, then $\lambda$ is the
cone generated by $v_1-z,\dots,v_k-z$; especially, $\lambda$
is a polyhedral cone, and it is easily seen that $\La\lambda\Ra=W$.
Moreover, $\lambda$ is strictly convex; indeed, otherwise there
exist real numbers $a_1,\dots,a_k\geq 0$, with $a_i>0$ for at
least an index $i\leq k$, such that $\sum_{i=1}^ka_i\cdot(v_i-z)=0$,
{\em i.e.} $\sum_{i=1}^ka_iv_i=(\sum^k_{i=1}a_i)\cdot z$, which
is absurd, since $z\notin\sigma$. Pick $u\in\lambda^\vee$ such
that $\lambda\cap\Ker\,u=0$, and set
$C_\lambda:=\lambda\cap u^{-1}(1)$; by lemma \ref{lem_cone-cell},
the subset $C_\lambda$ is homeomorphic to $\bar\B{}^{d-2}(1)$.
Lastly, let $\pi:W\setminus\Ker\,u\to u^{-1}(1)$ be the radial
projection (so $\pi(w)$ is the intersection point of $\R w$ with
$u^{-1}(1)$, for every $w\in W\setminus\Ker\,u$). It is easily
seen that $\pi$ maps $S'$ bijectively onto $C_\lambda$, so the
restriction of $\pi$ is a homeomorphism $S'\isom C_\lambda$, as
required. 
\end{proof}

\begin{lemma}\label{lem_hide-and-seek}
In the situation of \eqref{subsec_back-to-cohomology}, let
$\lambda\subset\sigma$ be any face. For every $l\in L_\R$ we have:
\begin{enumerate}
\item
The set of points of $\sigma$ that are visible from $l$
is a union of faces of $\sigma$.
\item
Suppose $l\in L$, and endow $R[P_\lambda]$ with its natural
$L$-grading. Then
$$
\gr_lR[P_\lambda]=\left\{\begin{array}{ll}
          Rl & \text{if\/ $\lambda\subset\sigma$ is not visible from $l$} \\
          0 & \text{otherwise}.
                        \end{array}\right.
$$
\end{enumerate}
\end{lemma}
\begin{proof} Let $\rho_1,\dots,\rho_n$ be a system of generators
of $\sigma^\vee$; we may assume that, for some integer $k\leq n$
we have
$$
\rho_i(l)<0
\qquad
\text{if and only if $1\leq i\leq k$}.
$$
Arguing as in the proof of lemma \ref{lem_do-you-see-me}, we
check easily that the set of points of $\sigma$ visible from
$l$ equals $S:=\sigma\cap\bigcup_{i=1}^k\Ker\,\rho_i$, which
already shows (i).

(ii): Suppose first that $\lambda\subset\sigma$ is not visible from
$l$; then there exists $x\in\lambda\setminus S$, and since $F_\lambda$
generates $\lambda$ (see \eqref{subsec_from-con-to-mon}), we
may assume that $x\in F_\lambda$ (details left to the reader).
A simple inspection then shows that there exists a sufficiently
large $N\in\N$ such that $l+Nx\in\sigma$, whence $l\in R[P_\lambda]$,
and so $\gr_lR[P_\lambda]=Rl$. Conversely, if the latter identity
holds, then there exists $x\in F_\lambda$ such that $l+x\in P$,
whence $\rho_i(x)>0$ for $i=1,\dots,k$, so $x$ is not visible from $l$.
\end{proof}

\sset\subsubsection{}\label{subsec_hide-and-seek}
In the situation of \eqref{subsec_back-to-cohomology},
let us fix a linear form $u_0\in\sigma^\vee$ such that
$\sigma\cap\Ker\,u_0=0$, and define $\sigma^\circ$ and
$C_\sigma$ as in \eqref{subsec_cone-cell}. For any $l\in L_\R$,
denote by $S_l$ the set of points of $\sigma$ that are visible
from $l$, so $S_l$ is a union of faces of $\sigma$, by lemma
\ref{lem_hide-and-seek}(i), and therefore $C_l:=S_l\cap C_\sigma$
is a subcomplex of $C^\bullet_\sigma$.

\begin{proposition}\label{prop_hide-and-seek}
With the notation of \eqref{subsec_hide-and-seek}, suppose
$l\in L$. Then the following holds :
\begin{enumerate}
\item
If $-l\in\sigma^\circ$, the complex of $R$-modules
$\gr_l\bar\cC{}^\bullet_P$ is isomorphic to $R[-d]$.
\item
If $-l\notin\sigma^\circ$, there is a natural isomorphism
of complexes of $R$-modules :
$$
\gr_l\bar\cC{}^\bullet_P\isom\Hom^\bullet_\Z
(\bar\cC_\bullet(C^\bullet_\sigma)/\bar\cC_\bullet(C^\bullet_l),R[-1]).
$$
Moreover, in this case, both $\bar\cC_\bullet(C^\bullet_\sigma)$
and $\bar\cC_\bullet(C^\bullet_l)$ are acyclic complexes.
\end{enumerate}
\end{proposition}
\begin{proof}(i): If $-l\in\sigma^\circ$, then $l\in P_\lambda$
if and only if $\lambda=\sigma$, so the assertion is clear.

(ii): The sought identification of complexes follows from
lemma \ref{lem_hide-and-seek}(ii), by a direct inspection
of the constructions (and indeed, this holds even if
$-l\in\sigma^\circ$ : details left to the reader). Moreover,
it is already known from proposition \ref{prop_computes-sing-hom}
and lemma \ref{lem_cone-cell} that $\bar\cC_\bullet(C^\bullet_\sigma)$
is acyclic.

$\bullet$\ \
Next, if $l\in P$, then clearly $S_l=\emptyset$, so the
assertion for $\bar\cC_\bullet(C^\bullet_l)$ is trivial
in this case.

$\bullet$\ \
Thus, suppose that $l\notin P$; if furthermore $-l\notin P$,
then the convex cone $\tau$ generated by $P$ and $l$ is still
strictly convex, so we may find $u_1\in\tau^\vee$ such that
$u_1(l)=1$ and $\sigma\cap\Ker\,u_1=0$. Set
$C'_\sigma:=\sigma\cap u_1^{-1}(1)$ and $C'_l:=S_l\cap C'_\sigma$.
In remark \ref{rem_cone-complex}(i) we have exhibited a homeomorphism
$C{}^{\prime\bullet}_\sigma\isom C_\sigma^\bullet$ that preserves
the respective cell complex structures, and a simple inspection
shows that this homeomorphism maps $C'_l$ onto $C_l$. On the
other hand, it is easily seen that $C'_l$ is also the set of
points of $C'_\sigma$ that are visible from $l$ (indeed, the
segment that joins any point of $C'_\sigma$ to $l$ lies in
$u_1^{-1}(1)$, so its intersection with $\sigma$ equals its
intersection with $C'_\sigma$ : details left to the reader).
By virtue of lemma \ref{lem_do-you-see-me}(ii) (and proposition
\ref{prop_computes-sing-hom}), it follows that
$\bar\cC_\bullet(C^{\prime\bullet}_l)$ is acyclic, so the
same holds for $\bar\cC_\bullet(C^\bullet_l)$.

$\bullet$\ \
Lastly, suppose $-l\in P\setminus\sigma^\circ$; we let
$\rho_1,\dots,\rho_n$ be a system of generators of $\sigma^\vee$,
and $k\leq n$ an integer such that $\rho_i(l)<0$ if and only if
$1\leq i\leq k$. Denote by $\tau^\vee$ (resp. $\mu^\vee$) the
convex cone in $L_\R^\vee$ generated by $(\rho_{k+1},\dots,\rho_n)$
(resp. by $(-\rho_1,\dots,-\rho_k)$), and let $\tau$ (resp. $\mu$)
be the dual of $\tau^\vee$ (resp. of $\mu^\vee$) in $L_\R$;
with this notation, we have $l\in\tau\cap\mu^\circ$.
Especially, $\tau\cap\mu^\circ\neq\emptyset$, and since
$\tau^\circ$ is dense in $\tau$ (by proposition
\ref{prop_spanning-cone}(i)), we deduce that
$\tau^\circ\cap\mu^\circ\neq\emptyset$ as well. Pick
$z\in\tau^\circ\cap\mu^\circ$; by inspecting the proof of
lemma \ref{lem_hide-and-seek}(i), it is easily seen that
$S_l=S_z$. But by construction, $-z\notin\sigma$, therefore
-- arguing as in the previous case -- we conclude that
$\bar\cC_\bullet(C^\bullet_z)=\bar\cC_\bullet(C^\bullet_l)$
is acyclic.
\end{proof}

\begin{theorem}[Hochster]\label{th_Hochster}
Let $R$ be a Cohen-Macaulay noetherian ring, $P$ a fine, sharp
and saturated monoid. Then $R[P]$ is a Cohen-Macaulay ring.
\end{theorem}
\begin{proof} In view of \cite[p.181, Cor.]{Mat} we may assume
that $R$ is a field. We argue by induction on $d:=\dim P$.
If $d=0$, we have $P=0$, and there is nothing to show. Suppose
that $d>0$ and that the assertion is already known for every
Cohen-Macaulay ring $R$ and every monoid as above, of dimension
$<d$. Set $L:=P^\gp$, and let $\sigma\subset L_\R$ be the unique
convex polyhedral cone such that $P=L\cap\sigma$. The ideal
$\fn:=R[\fm_P]$ is maximal in $R[P]$, and proposition
\ref{prop_hide-and-seek} and theorem
\ref{th_combinatorial-coh-compute} show that
\set\begin{equation}\label{eq_locally-CM}
\depth\,R[P]_\fn=\dim\,P.
\end{equation}
On the other hand, we have the more general :

\begin{claim}\label{cl_NullSt}
Let $F$ be a field, $P$ a fine monoid, $A$ an integral domain which
is an $F$-algebra of finite type, and $F'$ the field of fractions of $A$.
Then we have :
\begin{enumerate}
\item
For every maximal ideal $\fm\subset A$, the Krull dimension of
$A_\fm$ equals the transcendence degree of $F'$ over $F$.
\item
For every maximal ideal $\fm\subset F[P]$, the Krull dimension
of $F[P]_\fm$ equals $\rk_\Z\,P^\gp$.
\end{enumerate}
\end{claim}
\begin{pfclaim}(i): This is a straightforward consequence of
\cite[Th.5.6]{Mat}.

(ii): Choose an isomorphism : $P^\gp\isom L\oplus T$, where $T$ is
the torsion subgroup of $P^\gp$, and $L$ is a free abelian group of
finite rank; there follows an induced isomorphism
\eqref{eq_push-out-tensor} :
$$
F[P^\gp]\isom F[L]\otimes_F F[T].
$$
Let $B$ be the maximal reduced quotient of $F[P^\gp]$ (so the kernel
of the projection $F[P^\gp]\to B$ is the nilradical); we deduce
that $B$ is a direct product of the type $\prod_{i=1}^n F'_i[L]$,
where each $F'_i$ is a finite field extension of $F$. By (i), the
Krull dimension of $F'_i[L]_\fm$ equals $r:=\rk_\Z P^\gp$ for every
maximal ideal $\fm\subset F'_i[L]$, hence every irreducible component
of $\Spec\,B$ has dimension $r$. Let also $C$ be the maximal reduced
quotient of $F[P]$; the natural map $C\to B$ is an injective
localization, obtained by inverting a finite system of generators of
$P$, hence the induced morphism $\Spec\,B\to\Spec\,C$ is an open
immersion with dense image. Let $Z$ be any (reduced) irreducible
component of $\Spec\,C$; again by (i) it follows that every
non-empty open subset of $Z$ has dimension equal to $\dim Z$, so
necessarily the latter equals $r$.
\end{pfclaim}

From \eqref{eq_locally-CM}, corollary \ref{cor_consequent}(i)
and claim \ref{cl_NullSt}(ii) we see already that $R[P]_\fn$ is
a Cohen-Macaulay ring. Next, let $\lambda_1,\dots,\lambda_k$
be the one-dimensional faces of $\sigma$, and define $P_{\lambda_i}$
as in \eqref{subsec_back-to-cohomology}, for every $i=1,\dots,k$.
In light of claim \ref{cl_radical-of-I}, we have
$$
\Spec\,R[P]\setminus\{\fn\}=\bigcup_{i=1}^k\Spec\,R[P_{\lambda_i}]
$$
so it remains to check that $R[P_{\lambda_i}]$ is Cohen-Macaulay
for every $i=1,\dots,k$. However, we may find a decomposition
$P_{\lambda_i}\isom Q_i\times G_i$, where $G_i$ is a free abelian
group of finite type, and $Q_i$ is a fine, sharp and saturated
monoid of dimension $d-1$ (lemma \ref{lem_decomp-sats}). Set
$S_i:=R[G_i]$; then $R[P_{\lambda_i}]=S_i[Q_i]$, and $S_i$ is a
Cohen-Macaulay ring (\cite[Th.17.7]{Mat}), so the claim follows
by inductive assumption.
\end{proof}

\sset\subsubsection{}
In the situation of \eqref{subsec_back-to-cohomology},
we endow $R[P]$ with its natural $L$-grading, and denote
by $\underline A:=(R[P],\gr_\bullet R[P])$ the resulting
$L$-graded $R$-algebra. Let $(M,\gr_\bullet M)$ be an
$L$-graded $\underline A$-module (see definition
\ref{def_Gamma-graded-algs}(ii)); we set
$$
M^\dagger:=\bigoplus_{l\in L}\Hom_R(\gr_lM,R)
$$
which is naturally an $R$-submodule of $M^*:=\Hom_R(M,R)$ :
indeed, any $R$-linear map $\gr_lM\to R$ yields a linear form
$M\to R$, after composition with the projection $M\to\gr_lM$.
Notice that $M^*$ is naturally an $R[P]$-module; namely for any
linear form $f:M\to R$ and any $x\in P$, one defines
$x\cdot f:M\to R$ by the rule : $x\cdot f(m):=f(xm)$ for every
$m\in M$. Then, it is easily seen that $M^\dagger$ is an
$R[P]$-submodule of $M^*$; more precisely, we have
\set\begin{equation}\label{eq_stability-under-P}
x\cdot\Hom_R(\gr_lM,R)\subset\Hom_R(\gr_{l-x}M,R)
\qquad
\text{for every $x\in P$}.
\end{equation}
In light of \eqref{eq_stability-under-P}, it is convenient
to define an $L$-grading on $M^\dagger$ by the rule :
$$
\gr_lM^\dagger:=\Hom_R(\gr_{-l}M,R)
\qquad
\text{for every $l\in L$}
$$
and then $(M^\dagger,\gr_\bullet M^\dagger)$ is naturally
an $L$-graded $\underline A$-module.

\begin{remark}\label{rem_reflex-graded}
(i)\ \
Clearly, the rule $M\mapsto M^\dagger$ yields a functor
from the category $\underline A\Mod$ of $L$-graded
$\underline A$-modules to $\underline A\Mod^o$.

(ii)\ \
Let $\underline A\Mod_\mathrm{rflx}$ denote the full
subcategory of $\underline A\Mod$ whose objects are
all the $L$-graded $\underline A$-modules $(M,\gr_\bullet M)$
such that $\gr_lM$ is a projective $R$-module of finite type,
for every $l\in L$. A simple inspection of the definition
shows that the functor $(-)^\dagger$ of (i) restricts to an
equivalence of categories :
$$
\underline A\Mod_\mathrm{rflx}\to(\underline A\Mod_\mathrm{rflx})^o.
$$
\end{remark}

\begin{example}\label{ex_dagger-duals}
If $S\subset L$ is any $P$-submodule, notice that $L\setminus(-S)$
is also a $P$-submodule of $L$, and set
$$
S^\dagger:=L/(L\setminus(-S))
$$
where the quotient is a pointed $P$-module, as in remark
\ref{rem_point-and-complete}(iii). Explicitly, $S^\dagger$
is the set $(-S)\cup\{0_{S^\dagger}\}$ (where the zero element
$0_{S^\dagger}$ should not be confused with the neutral element
$0$ of the abelian group $L$); the $P$-module structure on
$S^\dagger$ is determined by the rule :
$$
x\cdot s:=\left\{\begin{array}{ll}
                 x+s        & \text{if $x+s\in -S$} \\
                 0_{S^\dagger} & \text{otherwise}
                 \end{array}\right.
\qquad
\text{for every $s\in -S$}.
$$
Then it is easily seen that there exists a natural isomorphism
of $L$-graded $\underline A$-modules :
$$
R[S]^\dagger\isom R\La S^\dagger\Ra
$$
(notation of \eqref{subsec_LaRa}). On the other hand, the
natural map $R[S]\to(R[S]^*)^*$ induces an isomorphism
$R[S]\isom (R[S]^\dagger)^\dagger$ whence an isomorphism of
$L$-graded $\underline A$-modules
\set\begin{equation}\label{eq_double-dagger}
R[S]\isom R\La S^\dagger\Ra^\dagger.
\end{equation}
\end{example}

\begin{theorem}[Stanley, Danilov]\label{th_Stanley}
In the situation of \eqref{subsec_back-to-cohomology}, set
$P^\circ:=L\cap\sigma^\circ$. We have :
\begin{enumerate}
\item
There exists a natural isomorphism :
$$
R\Hom^\bullet_{R[P]}(R,R[P^\circ])\isom R[-d]
\qquad
\text{in $\sD(R[P]\Mod)$}
$$
(here $R$ is regarded as an $R[P]$-module, via the augmentation
map $R[P]\to R$).
\item
If $R$ is a Gorenstein noetherian ring, the complex of coherent
$\cO_{\!X}$-modules $R[P^\circ]^\sim[0]$ is dualizing on $X$.
\end{enumerate}
\end{theorem}
\begin{proof}(i): From proposition \ref{prop_hide-and-seek} we
deduce a map of complexes of $L$-graded $\underline A$-modules
$$
\phi^\bullet:\bar\cC{}^\bullet_P\to R\La(P^\circ)^\dagger\Ra[-d]
$$
(notation of example \ref{ex_dagger-duals}) with
$\gr_l\phi^\bullet$ a quasi-isomorphism of complexes of
$R$-modules, for every $l\in L$. Since $\gr_l\bar\cC{}^\bullet_P$
is a finite dimensional $R$-vector space for every $l\in L$,
there follows -- in view of \eqref{eq_double-dagger} -- a
quasi-isomorphism of complexes of $L$-graded $\underline A$-modules
\set\begin{equation}\label{eq_complex-of-daggers}
R[P^\circ][0]\isom(\bar\cC{}^\bullet_P)^\dagger[-d].
\end{equation}
whence  a convergent spectral sequence
$$
E^{pq}_1:=\Ext^p_{R[P]}(R,(\bar\cC{}^{d-q}_P)^\dagger)\Rightarrow
\Ext^{p+q}_{R[P]}(R,R[P^\circ]).
$$
Hence, the assertion follows immediately from the :

\begin{claim}\label{cl_envelope-dagger}
With the foregoing notation, we have :
\begin{enumerate}
\item
$E^{pq}_1=0$ whenever $q<d$.
\item
$E^{0d}_1\simeq R$, and $E^{pd}_1=0$ for every $p>0$.
\end{enumerate}
\end{claim}
\begin{pfclaim}(i): Let $\lambda\subset\sigma$ be any face
with $\dim\lambda>0$; it suffices to check that
$E^i:=\Ext^i_{R(P]}(R,R[P_\lambda]^\dagger)=0$ for every $i\in\N$.
To this aim, pick any $x\in F_\lambda\setminus\{0\}$; since
$R[P_\lambda]^\dagger=R\La P^\dagger_\lambda\Ra$ (example
\ref{ex_dagger-duals}), we see that scalar multiplication
by $x$ is an automorphism on $R[P_\lambda]^\dagger$, hence
also on $E^i$. On the other hand, scalar multiplication by
$x$ is the zero endomorphism on $R$, hence also on $E^i$,
and the claim follows.

(ii): Pick a resolution $F_\bullet$ of the $P$-graded
$\Z[P]$-module $\Z:=\Z[P]/\Z[\fm_P]$ as in remark
\ref{rem_graded-resol}(iii), so that $F_n$ is a free
$\Z[P]$-module of finite type for every $n\in\N$, and
the differential $d_n:F_n\to F_{n-1}$ is a morphism of
$P$-graded $\Z[P]$-modules, for every $n>0$. Since the
augmented complex $F_\bullet\to\Z$ is flat as a complex
of $\Z$-modules, the complex of $R[P]$-modules
$R\otimes_\Z F_\bullet$ is still a free resolution of the
$R[P]$-module $R$, and we get natural isomorphisms
$$
E^{pd}_1\isom
H^p\Hom^\bullet_{R[P]}(R\otimes_\Z F_\bullet,R[P]^\dagger[0])
\isom H^p((R\otimes_\Z F_\bullet)^\dagger)
\qquad
\text{for every $p\in\N$}.
$$
(where the differentials of the complex
$(R\otimes_\Z F_\bullet)^\dagger$ are the maps $d_\bullet^\dagger$).
However, remark \ref{rem_reflex-graded}(ii) implies that
the complex $(R\otimes_\Z F_\bullet\to R)^\dagger$ is still
acyclic, whence the contention.
\end{pfclaim}

(ii): To ease notation, set $S:=\Spec\,R$ and
$\omega_P:=R[P^\circ]^\sim$; in view of proposition \ref{prop_fishy},
it suffices to check that the complex of coherent
$\cO_{\!X(x)}$-modules $\omega_P(x)[0]$ is dualizing on $X(x)$, for
every $x\in X$ (notation of definition \ref{def_strict-loc}(iii)).
Hence, fix any such point $x$, and let $\fp\subset R[P]$ be
the prime ideal corresponding to $x$; after replacing $R$
by its localization at the prime ideal $\fp\cap R$, we may
assume that $R$ is a local ring, and the structure morphism
$X\to S$ maps $x$ to the closed point of $S$, corresponding
to the maximal ideal $\fm_R$ of $R$. Denote by
$\lambda\subset\sigma$ the unique face such that
$P\setminus F_\lambda=\fp\cap P$, so that
$x\in U_\lambda:=\Spec\,R[P_\lambda]$. We may find a
decomposition $P\isom F_\lambda^\gp\times Q$, where $Q$ is also a
fine, sharp and saturated monoid (lemma \ref{lem_decomp-sats}),
whence an isomorphism of $S$-schemes
$$
U_\lambda\isom\Spec\,R[F^\gp_\lambda]\times_SY
\qquad
\text{where $Y:=\Spec\,R[Q]$}
$$
and by construction, $F^{-1}(\fp\cap P)$ is the maximal
ideal of $P_\lambda$, so the induced projection $p:U_\lambda\to Y$
maps $x$ to the maximal ideal $R[\fm_Q]+\fm_R[Q]\subset R[Q]$.
Let $\tau\subset Q^\gp_\R$ be the unique polyhedral cone such
that $Q=\tau\cap Q^\gp$, set $Q^\circ:=Q\cap\tau^\circ$, and
define the coherent $\cO_Y$-module $\omega_Q:=R[Q^\circ]^\sim$.
It is easily seen that there is a natural identification
$$
\omega_{P|U_\lambda}\isom p^*\omega_Q.
$$
In view of proposition \ref{prop_Goren-pullback} and remark
\ref{rem_Gorenstein}(i), it then suffices to check that
$\omega_Q(p(x))$ is dualizing on $Y(p(x))$. Thus, we may
replace $P$ by $Q$, and assume from start that
$\fp=R[\fm_P]+\fm_R[P]\subset R[P]$. In this case, set
$S_0:=\Spec\,R/\fm_R$, and denote by $i_1:S_0\to S$ and
$i_2:S\to X(S)$ the closed immersions (where $i_2$ corresponds
to the ring homomorphism $R(P]_\fp\to R$ deduced from
augmentation map); according to (i) and corollary
\ref{cor_Ext-loc=glob}(ii), we have a natural identification
$$
i_2^!\omega_P(x)[0]\simeq\cO_{\!S}[-d].
$$
On the other hand, since $R$ is Gorenstein, $\cO_{\!S}[-d]$
is a dualizing complex on $S$; from lemma \ref{lem_transit-dual}(i)
and proposition \ref{prop_sharp-flat}(i) it follows that
$i^!_1\cO_{\!S}[-d]\isom(i_2\circ i_1)^!\omega_P(x)[0]$
is dualizing on $S_0$, and therefore it is isomorphic
to $\cO_{\!S_0}[c]$ for some $c\in\Z$. Lastly, proposition
\ref{prop_Stanley-criterion}(i) and corollary \ref{cor_Ext-loc=glob}(ii)
now say that $\omega_P(x)[0]$ is dualizing on $X(x)$.
\end{proof}

\begin{remark} Suppose that $R$ is a field; then example
\ref{ex_SGA2-for-coherents} shows that $(\bar\cC{}^0_P)^\dagger$
is the injective hull of the $R[P]$-module $R$. In this case,
claim \ref{cl_envelope-dagger}(ii) is therefore a special
case of claim \ref{cl_hom-to-hull}. For a general $R$, the
$R[P]$-module $(\bar\cC{}^0_P)^\dagger$ shall not be injective,
but it can be viewed as a graded variant of the injective hull
construction.
\end{remark}

\section{Logarithmic geometry}\label{chap_logschemes}
\subsection{Log topoi}
Henceforth, {\em all topoi under consideration will be locally
ringed and with enough points, and all morphisms of topoi will
be morphisms of locally ringed topoi} (see
\eqref{subsec_topolo-on-sch}). The purpose of this restriction
is to insure that we obtain the right notions, when we specialize
to the case of schemes.

\begin{definition}
Let $T:=(T,\cO_{\!T})$ be a locally ringed topos.

\begin{enumerate}
\item
A {\em pre-log structure\/} on $T$ is the datum of a pair
$$
(\underline{M},\alpha)
$$
where $\underline{M}$ is a $T$-monoid, and
$\alpha:\underline{M}\to\cO_{\!T}$ is a morphism of
$T$-monoids, called the {\em structure map\/} of
$\underline M$, and where the monoid structure on
$\cO_{\!T}$ is induced by the multiplication law (hence,
by the multiplication in the ring $\cO_{\!T}(U)$, for
every object $U$ of $T$).
\item
A morphism $(\underline M,\alpha)\to(\underline N,\beta)$ of
pre-log structures on $T$, is a map
$\gamma:\underline M\to\underline N$ of $T$-monoids, such that
$\beta\circ\gamma=\alpha$. We denote by
$$
\bprelog_T
$$
the category of pre-log structures on $T$.
\item
A pre-log structure $(\underline M,\alpha)$ on $T$ is 
a {\em log structure\/} if $\alpha$ restricts to an
isomorphism:
$$
\alpha^{-1}\cO^\times_{\!T}\isom
\underline M^\times\isom\cO^\times_{\!T}.
$$
The datum of a locally ringed topos $(T,\cO_{\!T})$ and a log
structure on $T$ is also called, for short, a {\em log topos}.
We denote by
$$
\blog_T
$$
the full subcategory of $\bprelog_T$ consisting of all log
structures on $T$. When there is no danger of ambiguity, we
shall often omit mentioning explicitly the map $\alpha$, and
therefore only write $\underline M$ to denote a pre-log or a
log structure.
\end{enumerate}
\end{definition}

\sset\subsubsection{}
The category $\blog_T$ admits an initial object, namely
the log structure $(\cO^\times_{\!T},j)$, where $j$ is the
natural inclusion; this is called the {\em trivial log structure\/}.
$\blog_T$ admits a final object as well: this is
$(\cO_{\!T},\one_{\cO_{\!T}})$.
A morphism of locally ringed topoi $f:T\to S$ induces a pair
of adjoint functors :
\set\begin{equation}\label{eq_log-up-and-down}
f^*:\bprelog_S\to\bprelog_T\qquad f_*:\bprelog_T\to\bprelog_S.
\end{equation}
Namely, let $(\underline M,\alpha)$ (resp. $(\underline N,\beta)$)
be a pre-log structure on $S$ (resp. on $T$) and
$$
f^\flat:\cO_S\to f_*\cO_T
\qquad
f^\sharp:f^{-1}\cO_S\to\cO_T
$$
the natural morphisms (corresponding to each other under the
adjunction $(f^{-1},f_*)$ that defines the morphism $f$); then
$$
f^{-1}\underline M\xrightarrow{\ f^{-1}\alpha\ }f^{-1}\cO_S
\xrightarrow{\ f^\sharp\ }\cO_T
$$
defines $f^*(\underline M,\alpha)$, and $f_*(\underline N,\beta)$
is the pair $(\underline P,\gamma)$, where $\underline P$
is the fibre product in the cartesian diagram
$$
\xymatrix{
\underline P \ar[r]^-\gamma \ar[d] & \cO_S \ar[d]^{f^\flat} \\
f_*\underline N \ar[r]^-{f_*\beta} & f_*\cO_T.
}$$

\begin{lemma}\label{lem_check-iso}
Let $T$ be a topos with enough points, $\phi:\cM\to\cN$
a morphism of integral $T$-monoids inducing isomorphisms
$\cM^\times\isom\cN^\times$ and $\cM^\sharp\isom\cN^\sharp$.
Then $\phi$ is an isomorphism.
\end{lemma}
\begin{proof} This can be checked on the stalks, hence
we are reduced to the corresponding assertions for a morphism
$M\to N$ of monoids. Moreover, $M^\sharp$ is just the
set-theoretic quotient of $M$ by the translation action of
$M^\times$ (lemma \ref{lem_special-p-out}(iii)), so the assertion is
straightforward, and shall be left as an exercise for the reader.
\end{proof}

\begin{definition}
Let $\gamma:(\underline M,\alpha)\to(\underline N,\beta)$
be a morphism of pre-log structures on the locally ringed
topos $T$, and $\xi$ a $T$-point.
\begin{enumerate}
\item
We say that $(\underline M,\alpha)$ is {\em integral\/} (resp.
{\em saturated}) if $\underline M$ is an integral (resp. integral
and saturated) $T$-monoid.
\item
We say that $\gamma$ is {\em flat} (resp. {\em saturated})
{\em at the point $\xi$}, if
$\gamma_\xi:\underline M{}_\xi\to\underline N{}_\xi$ is a flat
morphism of monoids (resp. a saturated morphism of integral monoids)
(see remark \ref{rem_flatness}(vi)).
\item
We say that $\gamma$ is {\em flat} (resp. {\em saturated}),
if $\gamma$ is a flat morphism of $T$-monoids  (resp. a
saturated morphism of integral $T$-monoids) (see definition
\ref{def_satur-top}).
In view of proposition \ref{prop_finally-flat} (resp.
corollary \ref{cor_satur-check}), this is the same as saying
that $\gamma$  is flat (resp. saturated) at every $T$-point.
\end{enumerate}
\end{definition}

\sset\subsubsection{}\label{subsec_up_and_down-log}
The forgetful functor :
$$
\blog_T\to\bprelog_T \quad : \quad
\underline{M}\mapsto\underline{M}^\prelog
$$
admits a left adjoint :
$$
\bprelog_T\to\blog_T\quad : \quad
(\underline{M},\alpha)\mapsto(\underline{M},\alpha)^{\log}
$$
such that the resulting diagram :
\set\begin{equation}\label{eq_cokart}
{\diagram
\alpha^{-1}(\cO^\times_{\!T}) \ar[r] \ar[d] & \underline{M} \ar[d] \\
\cO^\times_{\!T} \ar[r] & \underline{M}^{\log}
\enddiagram}
\end{equation}
is cocartesian in the category of pre-log structures.
One calls $\underline{M}^{\log}$ the {\em log structure
associated\/} to $\underline{M}$. From this description,
it is easily seen that the counit of adjunction
\set\begin{equation}\label{eq_log-counit}
(\underline M^\prelog)^{\log}\to\underline M
\end{equation}
is an isomorphism for every log structure $\underline M$
on $T$.

Composing with the adjunction \eqref{eq_log-up-and-down},
we deduce a pair of adjoint functors :
$$
f^*:\blog_S\to\blog_T\qquad f_*:\blog_T\to\blog_S
$$
for any map of locally ringed topoi $f:T\to S$. Explicitly, if
$\underline N$ is any log structure on $S$, then $f^*\underline N$
is the push-out in the cocartesian diagram of $T$-monoids :
$$
\xymatrix{
f^*\cO_{\!S}^\times \ar[r] \ar[d] &
f^*(\underline N^\prelog) \ar[d] \\
\cO^\times_T \ar[r] & f^*\underline N
}$$
and on the other hand, if $\underline M$ is any log
structure on $T$, then it is easily seen that the
pre-log structure $f_*(\underline M^\prelog)$ is
actually a log structure (this can be checked on the
stalks over the points of $T$).
It follows easily that the induced map of $T$-monoids :
\set\begin{equation}\label{eq_general-iso}
f^*(\underline N^\sharp)\to(f^*\underline N)^\sharp
\end{equation}
is an isomorphism. As a consequence, for every point $\xi$ of $T$,
the natural map $\underline N_{f(\xi)}\to(f^*\underline N)_\xi$ is
a local morphism of monoids.

\begin{remark}\label{rem_colim-log}
The category $\blog_T$ admits arbitrary colimits : indeed,
since the counit \eqref{eq_log-counit} is an isomorphism,
it suffices to construct such colimits in the category of
pre-log structures, and then apply the functor
$(-)\mapsto(-)^{\log}$ which preserves colimits, since it
is a left adjoint.
\end{remark}

\sset\subsubsection{}\label{subsec_pointed-log-topos}
Let $(\underline M,\alpha)$ be a pre-log structure on $T$.
The morphism $\alpha$ extends to a unique morphism of pointed
$T$-monoids $\alpha_\circ:\underline M{}_\circ\to\cO_{\!T}$,
whence a new pre-log structure
$$
(\underline M,\alpha)_\circ:=(\underline M{}_\circ,\alpha_\circ).
$$
Clearly, $(\underline M,\alpha)$ is a log structure if and
only if the same holds for $(\underline M,\alpha)_\circ$.
More precisely, for any pre-log structure $\underline M$
there is a natural isomorphism of log structures :
$$
(\underline M{}_\circ)^{\log}\isom(\underline M^{\log})_\circ.
$$
Furthermore, for any morphism $f:T\to S$ of topoi, we have
natural isomorphisms of pre-log (resp. log) structures
$$
f^*(\underline N{}_\circ)\isom(f^*\underline N)_\circ
\qquad
f_*(\underline M{}_\circ)\isom(f_*\underline M)_\circ
$$
for every pre-log (resp. log) structure $\underline N$ on
$S$ and $\underline M$ on $T$ (details left to the reader).

\begin{example}\label{ex_push-trivial}
(i)\ \ 
Let $T\to S$ be a morphism of topoi. Since the initial object
of a category is the empty coproduct, it follows formally that the
inverse image $f^*(\cO^\times_{\!S},j)$ of the trivial log structure
on $S$, is the trivial log structure on $T$.

(ii)\ \ 
Let $T$ be a topos, and $j_U:T\!/U\to T$ an open subtopos
(see example \ref{ex_localization-topos}(i)). Consider the subsheaf of
monoids $\underline{M}\subset\cO_T$ such that :
$$
\underline{M}(V):=
\{s\in\cO_T(V)~|~ s_{|U\times V}\in\cO^\times_T(U\times V)\}
\qquad
\text{for every object $V$ of
$T$}.
$$
Then it is easily seen that the natural map $\underline{M}\to\cO_T$
is a log structure on $T$. This log structure is (naturally
isomorphic to) the extension $j_{U*}\cO^\times_{\!U}$ of the trivial
log structure on $T\!/U$.

(iii)\ \
Let $U$ be any object of the topos $T$, and $\underline M$
a log structure on $T$. Since $\cO_{T\!/U}=(\cO_T)_{|U}$, it is
easily seen that the natural morphism of pre-log structures
$$
(\underline M^\prelog)_{|U}\to(\underline M{}_{|U})^\prelog
$$
is an isomorphism.

(iv)\ \
Let $\beta:\underline M\to\cO_T$ be a pre-log structure on a
topos $T$. Then $\beta^{-1}(0)\subset\underline M$ is an ideal,
and $\beta$ factors uniquely through the natural map
$\underline M\to\underline M/\beta^{-1}(0)$, and a pre-log
structure
$$
(\underline M,\beta)_\red:=(\underline M/\beta^{-1}(0),\bar\beta)
$$
called the {\em reduced pre-log structure\/} associated to
$\underline M$. As usual, we shall often write just
$\underline M{}_\red$ instead of $(\underline M,\beta)_\red$.
We say that $\beta$ is {\em reduced\/} if the induced morphism
of pre-log structures $\underline M\to\underline M{}_\red$ is an
isomorphism.

Suppose now that $\underline M$ is a log structure; then it
is easily seen that the same holds for $\underline M{}_\red$.
More precisely, since the tensor product is right exact (see
\eqref{eq_finally-flatnes}), for any pre-log structure
$\underline M$ the natural morphism of log structures
$$
(\underline M{}_\red)^{\log}\to(\underline M^{\log})_\red
$$
is an isomorphism.
\end{example}

\begin{lemma}\label{lem_integr-flat}
Let $\gamma:(\underline M,\alpha)\to(\underline N,\beta)$ be a
morphism of pre-log structures on $T$. We have:
\begin{enumerate}
\item
If $\underline M$ is integral (resp. saturated), then the same
holds for $\underline M^{\log}$.
\item
The unit of adjunction $\underline M\to\underline M^{\log}$
is a flat morphism.
\item
If $\gamma$ is flat (resp. saturated) at a $T$-point $\xi$,
the same holds for the induced morphism
$\gamma^{\log}:\underline M^{\log}\to\underline N^{\log}$ of log structures.
\item
Especially, if $\gamma$ is flat (resp. saturated), the same
holds for $\gamma^{\log}$.
\end{enumerate}
\end{lemma}
\begin{proof} In view of lemma \ref{lem_concern}(ii)
and proposition \ref{prop_finally-flat}, both (i) and (ii)
can be checked on stalks. Taking into account lemma
\ref{lem_concern}(i), we are reduced to showing the following.
Let $P$ be a monoid, $A$ a ring, $\beta:P\to(A,\cdot)$
a morphism of monoids; then the natural map
$P\to P':=P\otimes_{\beta^{-1}A^\times}A^\times$ is flat,
and if $P$ is integral (resp. saturated), the same holds
for $P'$. The first assertion follows easily from example
\ref{ex_flat-quots}(vi), and the second follows from remark
\ref{rem_why-not}(i) (resp. from corollary \ref{cor_exact}(ix)
and proposition \ref{prop_motivate-def}).

(iii): The map $\gamma^{\log}$ can be factored as the composition
of 
$$
\gamma':=\gamma\otimes_{\beta^{-1}\cO_T^\times}\cO^\times_T:
\underline N^{\log}\to\underline P:=
\underline M\otimes_{\beta^{-1}\cO_T^\times}\cO^\times_T
$$
and the natural unit of adjunction
$\gamma'':\underline P\to\underline P^{\log}=\underline M^{\log}$.
If $\gamma_\xi$ is flat, the same clearly holds for $\gamma'_\xi$,
and (ii) says that $\gamma''$ is flat, hence $\gamma^{\log}_\xi$
is flat in this case. Lastly, suppose that $\gamma_\xi$ is saturated,
and we wish to show that $\gamma^{\log}_\xi$ is saturated. Set
$P:=\alpha^{-1}_\xi\cO^\times_{T,\xi}$ and
$Q:=\beta^{-1}_\xi\cO^\times_{T,\xi}$. Then the induced map
$(P^{-1}\underline M{}_\xi)^\sharp\to
(Q^{-1}\underline N{}_\xi)^\sharp$ is saturated (lemma
\ref{lem_little}(ii,iii)). But the latter is the same as the
morphism $(\gamma_\xi^{\log})^\sharp$, and then also $\gamma_\xi$
is saturated, again by lemma \ref{lem_little}(iii).
\end{proof}

\begin{lemma}\label{lem_give-it-a-name}
Let $f:T'\to T$ be a morphism of topoi, $\xi$ a $T'$-point,
and $\gamma:(\underline M,\alpha)\to(\underline N,\beta)$ a
morphism of integral log structures on $T$.
The following conditions are equivalent :
\begin{enumerate}
\alphaenu
\item
$\gamma$ is flat (resp. saturated) at the $T$-point
$f(\xi)$.
\item
$f^*\gamma$ is flat (resp. saturated) at the $T'$-point $\xi$.
\item
$\gamma^\sharp_\xi$ is a flat (resp. saturated) morphism
of monoids.
\end{enumerate}
\end{lemma}
\begin{proof} The equivalence of (a) and (c) follows
from corollary \ref{cor_yet-another-flat} (resp. lemma
\ref{lem_little}(ii)). By the same token, (b) holds if
and only if $(f^*\gamma)^\sharp_\xi$ is flat (resp. saturated);
in light of the isomorphism \eqref{eq_general-iso}, the latter
condition is equivalent to (c).
\end{proof}

\sset\subsubsection{}\label{subsec_Konstant}
For any locally ringed topos $T$, let us write the objects
of $\Mnd/\Gamma(T,\cO_{\!T})$ in the form $(M,\phi)$, where
$M$ is any monoid, and $\phi:M\to\Gamma(T,\cO_{\!T})$ a morphism
of monoids. There is an obvious global sections functor :
$$
\Gamma(T,-):\bprelog_T\to\Mnd/\Gamma(T,\cO_{\!T})\quad :\quad
(\underline{N},\alpha)\mapsto
(\Gamma(T,\underline{N}),\Gamma(T,\alpha))
$$
which admits a left adjoint :
$$
\Mnd/\Gamma(T,\cO_{\!T})\to\bprelog_T\quad :\quad
(M,\phi)\mapsto(M,\phi)_T:=(M_T,\phi_T).
$$
Indeed, $M_T$ is the constant sheaf on $(T,C_T)$ with value $M$, and
$\phi_T$ is the composition of the map of constant sheaves
$M_T\to\Gamma(T,\cO_{\!T})_T$ induced by $\phi$, with the natural
map $\Gamma(T,\cO_{\!T})_T\to\cO_{\!T}$. Again, we shall often just
write $M_T$ to denote this pre-log structure.

After taking associated log structures, we deduce a left
adjoint :
\set\begin{equation}\label{eq_another-left-adj}
\Mnd/\Gamma(T,\cO_{\!T})\to\blog_T\quad :\quad (M,\phi)\mapsto
M_T^{\log}:=(M,\phi)_T^{\log}
\end{equation}
to the global sections functor. $M^{\log}_T$ is called the
{\em constant log structure\/} associated to $(M,\phi)$.

\begin{definition}\label{def_chart}
Let $T$ be a locally ringed topos, $(\underline{M},\alpha)$ a log
structure on $T$.
\begin{enumerate}
\item
A {\em chart for $\underline{M}$\/} is an object $(P,\beta)$
of $\Mnd/\Gamma(T,\cO_{\!T})$, together with a map of pre-log
structures $\omega_P:(P,\beta)_T\to\underline{M}$, inducing
an isomorphism on the associated log structures.
(Notation of \eqref{subsec_Konstant}.)
\item
We say that a chart $(P,\beta)$ is {\em finite} (resp. {\em
integral}, resp. {\em fine}, resp. {\em saturated}) if $P$ is a
finitely generated (resp. integral, resp. fine, resp. integral and
saturated) monoid.
\item
Let $\phi:\underline M\to\underline N$ be a morphism of log
structures on $T$. A {\em chart for $\phi$\/} is the datum of charts :
$$
\omega_P:(P,\beta)_T\to\underline M
\quad\text{and}\quad
\omega_Q:(Q,\gamma)_T\to\underline N
$$
for $\underline M$, respectively $\underline N$, and
a morphism of monoids $\theta:Q\to P$, fitting into
a commutative diagram :
$$
\xymatrix{
Q_T \ar[rr]^-{\theta^{\log}_T}
\ar[d]_{\omega_Q} & &
P_T^{\log} \ar[d]^{\omega_P} \\
\underline N \ar[rr]^\phi & & \underline M.
}$$
We say that such a chart is {\em finite\/} (resp. {\em integral},
resp. {\em fine}, resp. {\em saturated}) if the monoids $P$ and $Q$
are finitely generated (resp. integral, resp. fine, resp. integral
and saturated). We say that the chart is {\em flat} (resp.
{\em saturated}), if $\theta$ is a flat morphism of monoids
(resp. a saturated morphism of integral monoids).
\item
We say that $\underline M$ is {\em quasi-coherent\/}
(resp. {\em coherent\/}) if there exist a covering family
$(U_\lambda\to 1_T~|~\lambda\in\Lambda)$ of the final object
$1_T$ in $(T,C_T)$, and for every $\lambda\in\Lambda$, a chart
(resp. a finite chart) $(P_\lambda,\beta_\lambda)$ for
$\underline M{}_{|U_\lambda}$.
\item
We say that $(\underline{M},\alpha)$ is {\em quasi-fine} (resp.
{\em fine}) if it is integral and quasi-coherent (resp. and coherent).
\item
Let $\xi$ be any $T$-point. We say that a chart $(P,\beta)$ is
{\em local} (resp. {\em sharp}) {\em at the point $\xi$}, if the
morphism $\beta_\xi:P\to\cO_{T,\xi}$ is local (resp. if $P$
is sharp and $\beta_\xi$ is local).
\end{enumerate}
\end{definition}

\begin{lemma}\label{lem_simple-charts-top}
Let $f:T\to S$ be a morphism of locally ringed topoi, $\underline Q$
a log structure on $S$, and $\xi$ any point of $S$. The following
holds :
\begin{enumerate}
\item
If $\underline Q$ is quasi-coherent (resp. coherent, resp.
integral, resp. saturated, resp. quasi-fine, resp. fine),
then the same holds for $f^*\underline Q$.
\item
Suppose that $\underline Q$ is an integral log structure.
Then $\underline Q^\sharp$ is an integral $S$-monoid, and
$\underline Q$ is saturated if and only if\/
$\underline Q^\sharp$ is a saturated $S$-monoid.
\item
Suppose that $\underline{Q}$ is quasi-coherent. Then $\underline Q$
is integral (resp. integral and saturated, resp. fine, resp. fine
and saturated) if and only if there exist a covering family
$(U_\lambda\to 1_S~|~\lambda\in\Lambda)$ of the final object
of $S$, and for every $\lambda\in\Lambda$, an integral (resp.
integral and saturated, resp. fine, resp. fine and saturated) chart
$(P_\lambda)_{U_\lambda}\to\underline Q{}_{|U_\lambda}$.
\item
If $\underline{P}$ is any coherent log structure on $S$, and
$\omega:\underline P{}_\xi\to\underline Q{}_\xi$ is a map of
monoids, then :
\begin{enumerate}
\item
There exist a neighborhood $U$ of $\xi$ and a morphism
$\vartheta:\underline P{}_{|U}\isom\underline Q{}_{|U}$ of log
structures, such that $\vartheta_\xi=\omega$.
\item
Moreover, for any two morphisms
$\vartheta,\vartheta':\underline P{}_{|U}\isom\underline Q{}_{|U}$
with the property of {\em (a)}, we may find a smaller neighborhood
$V\to U$ of\/ $\xi$ such that $\vartheta_{|V}=\vartheta'_{|V}$.
\item
Especially, if $\underline Q$ is also coherent, and $\omega$ is an
isomorphism, we may find $\theta$ and a small enough $U$ as in
{\em (a)}, such that $\vartheta$ is an isomorphism.
\end{enumerate}
\item
If $M$ is a finitely generated monoid, and $\omega:M\to\cO_{\!S,\xi}$
a morphism of monoids, then we may find a neighborhood $U$ of\/
$\xi$ and a morphism $\theta:M_U\to\cO_{\!U}$ of\/ $S/U$-monoids,
such that $\theta_\xi=\omega$.
\end{enumerate}
\end{lemma}
\begin{proof}(i): If $\underline{Q}$ is the constant log structure
associated to a map of monoids $\alpha:Q\to\Gamma(S,\cO_{\!S})$,
then $f^*\underline{Q}$ is the constant log structure associated to
$\Gamma(S,f^\natural)\circ\alpha$ (where $f^\natural:\cO_{\!S}\to
f_*\cO_{\!T}$ is the natural map). The assertions concerning
quasi-coherent or coherent log structures are a straightforward
consequence. Next, suppose that $\underline{Q}$ is integral (resp.
saturated); we wish to show that $f^*\underline{Q}$ is integral
(resp. saturated). To this aim, let
$\underline{M}:=f^*(\underline{Q}^\prelog)$; by lemma
\ref{lem_T-satura}(i), $\underline{M}$ is an integral (resp.
saturated) $T$-monoid; then the assertion follows from lemma
\ref{lem_integr-flat}(i).

(ii): By lemma \ref{lem_concern}(ii) the assertion can be checked
on stalks. Hence, suppose that $\underline Q$ is integral; then
$\underline Q{}_\xi$ is integral by {\em loc.cit.}, consequently, the
same holds for $\underline Q{}_\xi/\underline Q{}^\times_\xi$ (lemma
\ref{lem_integral-quot}). The second assertion follows from lemma
\ref{lem_exc-satura}(ii).

(iii): Suppose first that $\underline{Q}$ is quasi-coherent and
integral. Hence, there is a covering family $(U_\lambda\to
1_S~|~\lambda\in\Lambda)$, and for every $\lambda\in\Lambda$ a
monoid $M_\lambda$, a pre-log structure
$\alpha_\lambda:(M_\lambda)_{U_\lambda}\to\cO_{\!U_\lambda}$, and
an isomorphism $((M_\lambda)_{U_\lambda},\alpha_\lambda)^{\log}
\isom\underline{Q}_{|U_\lambda}$; whence a cocartesian diagram
of $S$-monoids, as in \eqref{eq_cokart} :
\set\begin{equation}\label{eq_cokart-again}
{\diagram
\underline{N}:=\alpha^{-1}_\lambda(\cO^\times_{\!U_\lambda})
          \ar[r] \ar[d] & (M_\lambda)_{\!U_\lambda} \ar[d] \\
          \cO^\times_{\!U_\lambda} \ar[r] & \underline{Q}_{|U_\lambda}.
\enddiagram}
\end{equation}
The induced diagram $\eqref{eq_cokart-again}^\intg$ of integral
$S$-monoids is still cocartesian. Moreover, since
$\underline Q_{|U_\lambda}$ is integral, $\alpha_\lambda$ factors
through a unique map
$\beta_\lambda:((M_\lambda)_{U_\lambda})^\intg\to
\underline Q_{|U_\lambda}\to\cO_{\!U_\lambda}$,
and the morphism in $S$ underlying the induced morphism
of $S$-monoids $\underline{N}^\intg\to\underline{N}':=
\beta^{-1}_\lambda(\cO^\times_{\!U_\lambda})$
is an epimorphism (this can be checked easily on the stalks).
Furthermore,
$((M_\lambda)_{U_\lambda})^\intg\simeq(M_\lambda^\intg)_{U_\lambda}$
(see \eqref{eq_nat-int-glob}). It follows that the natural map
$$
\underline{Q}_{|U_\lambda}\to
\cO^\times_{\!U_\lambda}\amalg_{\underline{N}'}(M_\lambda^\intg)_{U_\lambda}
\simeq(M_\lambda^\intg)^{\log}_{U_\lambda}
$$
is an isomorphism, so the claim holds with
$P_\lambda:=M_\lambda^\intg$. If $Q$ is fine, we can find
$M_\lambda$ as above which is also finitely generated, in which case
the resulting $P_\lambda$ shall be fine.

Suppose additionally, that $\underline{Q}$ is saturated.
By the previous case, we may then find a covering family
$(U_\lambda\to e_S~|~\lambda\in\Lambda)$, and for every
$\lambda\in\Lambda$ an integral monoid $M_\lambda$, a pre-log
structure
$\alpha_\lambda:(M_\lambda)_{U_\lambda}\to\cO_{\!U_\lambda}$, and an
isomorphism $((M_\lambda)_{U_\lambda},\alpha_\lambda)^{\log}
\isom\underline{Q}_{|U_\lambda}$; whence a cocartesian diagram
\eqref{eq_cokart-again} consisting of integral $S$-monoids.
The induced diagram $\eqref{eq_cokart-again}^\sat$ is still
cocartesian; one may then argue as in the foregoing, to obtain a
natural isomorphism :
$\underline Q{}_{|U_\lambda}\isom(M_\lambda^\sat)^{\log}_{U_\lambda}$.
Furthermore, if $M_\lambda$ is finitely generated, the same holds
for $M_\lambda^\sat$ (corollary \ref{cor_fragment-Gordon}(ii)),
hence the chart thus obtained shall be fine and saturated, in
this case.

Conversely, if a family $(P_\lambda~|~\lambda\in\Lambda)$ of
integral (resp. saturated) monoids can be found fulfilling the
condition of (iii), then $(P_{\lambda})_{U_\lambda}$ is an integral
(resp. saturated) pre-log structure on $T\!/U_\lambda$ (example
\ref{ex_constant-mond}(ii)), hence the same holds for
its associated log structure $\underline Q{}_{|U_\lambda}$
(lemma \ref{lem_integr-flat}(i)), and thus also for
$\underline Q$ (lemma \ref{lem_concern}(ii), example
\ref{ex_covers-are-conserv}, and example
\ref{ex_push-trivial}(iii)). Moreover, if each $P_\lambda$
is fine, then $\underline Q$ is fine as well.

(iv.a): We may assume that $\underline P$ admits a finite chart
$\alpha:M_S\to\underline P$, for some finitely generated monoid
$M$, denote by $\beta:\underline Q\to\cO_{\!S}$ the structure map of
$\underline Q$, and set
$\omega':=\omega\circ\alpha_\xi:M\to\underline Q{}_\xi$. According to
lemma \ref{lem_finite-pres}(ii), the morphism $\omega'$ factors
through a map $\omega'':M\to\Gamma(U',\underline Q)$, for some
neighborhood $U'$ of $\xi$. By adjunction, $\omega''$ determines a
morphism of $S/\!U'$-monoids $\psi:M_{U'}\to\underline Q{}_{|U'}$
whence a morphism of pre-log structures
\set\begin{equation}\label{eq_sought-map}
(M_{U'},\beta_{|U'}\circ\psi)\to(\underline Q{}_{|U'},\beta_{|U'}).
\end{equation}
Let us make the following general observation :

\begin{claim}\label{cl_obviety}
Let $N$ be a finite monoid, $F$ any $S$-monoid, $f,g:N_S\to F$ two
morphisms of $S$-monoids, such that $f_\xi=g_\xi$. Then there exists
a neighborhood $U$ of $\xi$ in $S$ such that $f_{|U}=g_{|U}$.
\end{claim}
\begin{pfclaim} By adjunction, the morphisms $f$ and $g$ correspond
to unique maps of monoids $\Gamma(f),\Gamma(g):N\to\Gamma(S,F)$;
since $N$ is finite and $f_\xi=g_\xi$, we may find a neighborhood
$U$ of $\xi$ such that the maps $N\to F(U)$ induced by $\Gamma(f)$
and $\Gamma(g)$ coincide. Again by adjunction, we deduce a unique
morphism of $S/\!U$-monoids $N_U\to F_{|U}$, which by construction
is just the restriction of both $f$ and $g$.
\end{pfclaim}

Let $\gamma:\underline P\to\cO_{\!S}$ be the structure map of
$\underline P$; we apply claim \ref{cl_obviety} with $S$ replaced by
$S/\!U'$, to deduce that there exists a small enough neighborhood
$U\to U'$ of $\xi$ such that the restriction
$(\gamma\circ\alpha)_{|U}:M_{|U}\to\cO_{\!S|U}$ agrees with
$\beta_{|U}\circ\psi_{|U}$. Then it is clear that the morphism of
log structure associated to $\eqref{eq_sought-map}_{|U}$ yields the
sought extension $\theta$ of $\omega$.

(iv.b): By assumption we have the identity :
$\vartheta_\xi=\vartheta'_\xi$; however, any morphism of log
structures $\underline{P}_{|U}\to\underline{Q}_{|U}$ is already
determined by its restriction to the image of any finite local chart
$M_U\to\underline{P}_{|U}$. Hence the assertion follows from claim
\ref{cl_obviety}.

(iv.c): We apply (iv.a) to $\omega^{-1}$ to deduce the
existence of a morphism
$\sigma:\underline Q{}_{|U}\to\underline P{}_{|U}$ such that
$\sigma_\xi=\omega^{-1}$ on some neighborhood $U$ of $\xi$.
Hence, $(\vartheta\circ\sigma)_\xi=\one_{\underline Q{}_\xi}$
and $(\sigma\circ\vartheta)_\xi=\one_{\underline P{}_\xi}$.
By (iv.b), these identities persist on some smaller neighborhood.

(v): The proof is similar to that of (iv.a), though simpler :
we leave it as an exercise for the reader.
\end{proof}

\begin{definition}\label{def_amorph-log}
(i)\ \
A {\em morphism $(T,\underline M)\to(S,\underline N)$ of
topoi with pre-log} (resp. {\em log}) {\em structures},
is a pair $f:=(f,\log f)$ consisting of a morphism of locally
ringed topoi $f:T\to S$, and a morphism
$$
\log f:f^*\underline{N}\to\underline{M}
$$
of pre-log structures (resp. log structures) on $T$. We say that
$f$ is {\em log flat} (resp. {\em saturated}) if $\log f$ is a
flat (resp. saturated) morphism of pre-log structures.

(ii)\ \
Let $(f,\log f)$ as in (i) be a morphism of log topoi, $\xi$
a $T$-point; we say that $f$ is {\em strict at the point $\xi$},
if $\log f_\xi$ is an isomorphism. We say that $f$ is {\em strict},
if it is strict at every $T$-point.

(iii)\ \
A {\em chart\/} for $\phi$ is the datum of charts
$$
\omega_P:(P,\beta)_T\to\underline M
\quad\text{and}\quad
\omega_Q:(Q,\gamma)_S\to\underline N
$$
for $\underline M$, respectively $\underline N$, and a morphism
of monoids $\vartheta:Q\to P$, such that
$(f^*\omega_Q,\omega_P,\theta)$ is a chart for the morphism
$\log f$. We say that such a chart $(\omega_Q,\omega_P,\theta)$
is {\em finite} (resp. {\em fine}, resp. {\em flat}, resp.
{\em saturated}) if the same holds for the corresponding chart
$(f^*\omega_Q,\omega_P,\theta)$ of $\log f$.
\end{definition}

\begin{remark}\label{rem_strict-locus}
(i)\ \
Let $f:(T,\underline M)\to(S,\underline N)$ be a morphism
of log topoi, $g:T'\to T$ a morphism of topoi, and
$f':(T',g^*\underline M)\to(S,\underline N)$ the composition
of $f$ and the natural morphism of log topoi
$(T',g^*\underline M)\to(T,\underline M)$; let also $\xi$
be a $T'$-point. Then $f'$ is strict at the point $\xi$ if
and only if $f$ is strict at the point $g(\xi)$. Indeed,
$f'$ is strict at $\xi$ if and only if $(\log f')^\sharp_\xi$
is an isomorphism (lemma \ref{lem_check-iso}), if and only if
$(\log f)^\sharp_{g(\xi)}$ is an isomorphism (by
\eqref{eq_general-iso}), if and only if $\log f_{g(\xi)}$ is
an isomorphism (again by lemma \ref{lem_check-iso}).

(ii)\ \
For any log topos $(T,\underline M)$, let us set
$(T,\underline M)_\circ:=(T,\underline M{}_\circ)$.
Then $(T,\underline M)_\circ$ is a log topos (see
\eqref{subsec_pointed-log-topos}), and clearly, every
morphism $f:(T,\underline M)\to(S,\underline N)$ of log
topoi extends naturally to a morphism
$f_\circ:(T,\underline M)_\circ\to(S,\underline N)_\circ$
of log topoi.
\end{remark}

\sset\subsubsection{}\label{subsec_morph-logtop}
The $2$-category of log topoi admits arbitrary $2$-limits
indexed by usual categories. Indeed, if
$\cT:=((T_\lambda,\underline M{}_\lambda)~|~\lambda\in\Lambda)$
is any pseudo-functor (from a small category $\Lambda$, to
the category of log topoi), the $2$-limit of $\cT$ is the pair
$(T,\underline M)$, where $T$ is the $2$-limit of the system
of ringed topoi $(T_\lambda~|~\lambda\in\Lambda)$, and
$\underline M$ is the colimit of the induced system of log
structures $(p_\lambda^*\underline M_\lambda~|~\lambda\in\Lambda)$
on $T$, where $p_\lambda:T\to T_\lambda$ denotes the natural
projection, for every $\lambda\in\Lambda$ (see remark
\ref{rem_colim-log}).

\sset\subsubsection{}\label{subsec_ind-by-log-g}
Consider a $2$-cartesian diagram of log topoi :
$$
\xymatrix{(T',\underline M') \ar[r]^g \ar[d]_{f'} &
(T,\underline M) \ar[d]^f \\
(S',\underline N') \ar[r]^h & (S,\underline N).
}$$
The following result yields a relative variant of the
isomorphism \eqref{eq_general-iso} :

\begin{lemma}\label{lem_relative-var}
In the situation of \eqref{subsec_ind-by-log-g}, the morphism
$$
g^*\Coker(\log f)\to\Coker(\log f')
$$
induced by $\log g$ is an isomorphism of $T'$-monoids.
\end{lemma}
\begin{proof} Indeed, denote by $\beta:\underline M\to\cO_T$
and $\gamma:\underline N'\to\cO_{\!S'}$ the log structures
of $T$ and $S'$. Fix any $T'$-point $\xi$, let $\xi:=g(\xi')$,
and set
$$
P:=\underline M{}_\xi\otimes_{\underline N{}_{f(\xi)}}
\underline N'{}_{f'(\xi')}
\qquad\text{and}\qquad
\rho:=\beta_\xi\otimes\gamma_{f'(\xi')}:P\to\cO_{T'}.
$$
Then $\underline M'{}_{\xi'}=
P\otimes_{\rho^{-1}\cO_{T',\xi'}}\cO^\times_{T',\xi'}$, and
it is easily seen that
$\rho^{-1}\cO_{T',\xi'}=
\cO^\times_{T,\xi}\otimes_{\cO^\times_{\!S,f(\xi)}}
\cO^\times_{\!S',f'(\xi')}$. Therefore
$(\underline M'_{\xi'})^\sharp=P/\rho^{-1}\cO_{T',\xi'}=
\underline M{}^\sharp_\xi\otimes_{\underline N{}^\sharp_{f(\xi)}}
\underline N'{}^\sharp_{f'(\xi')}$, and
$$
\Coker(\log f'_{\xi'})=
\Coker(\underline N'{}^\sharp_{f'(\xi')}\to
\underline M{}^\sharp_\xi
\otimes_{\underline N{}^\sharp_{f(\xi)}}
\underline N'{}^\sharp_{f'(\xi')})=
\Coker(\underline N{}^\sharp_{f(\xi)}\to
\underline M{}^\sharp_\xi)
$$
whence the contention.
\end{proof}

\sset\subsubsection{}\label{subsec_chart-later} We consider
now a special situation that we will encounter in proposition
\ref{prop_good-charts}. Namely, let $Q$ be a monoid, and
$H\subset Q^\times$ a subgroup. Let also $G$ be an abelian
group, $\rho:G\to Q^\gp$ a group homomorphism, and set :
$$
H_\rho:=G\times_{Q^\gp}H \qquad Q_\rho:=G\times_{Q^\gp}Q.
$$
The natural inclusion $H\to Q$ and the projection $Q_\rho\to Q$
determine a unique morphism :
\set\begin{equation}\label{eq_chart-later}
Q_\rho\otimes_{H_\rho}H\to Q.
\end{equation}

\begin{lemma}\label{lem_charts-later}
In the situation of \eqref{subsec_chart-later}, suppose furthermore
that the composition :
$$
G\xrightarrow{\ \rho\ }Q^\gp\to(Q/H)^\gp
$$
is surjective. Then \eqref{eq_chart-later} is an isomorphism.
\end{lemma}
\begin{proof} Set $G':=G\oplus H$, and let $\rho':G'\to Q^\gp$
be the unique group homomorphism that extends $\rho$ and the
natural map $H\to Q\to Q^\gp$. Under the standing assumptions,
$\rho'$ is clearly surjective. Define $Q_{\rho'}$ and $H_{\rho'}$
as in \eqref{subsec_chart-later}; there is a natural isomorphism
of monoids : $Q_{\rho'}\isom Q_\rho\oplus H$, inducing an isomorphism
$H_{\rho'}\isom H_\rho\oplus H$, and defined as follows. To
every $g\in G$, $h\in H$, $q\in Q$ such that $[(g,h),q]\in Q_{\rho'}$,
we assign the element $[(g,h^{-1}q),a]\in Q_\rho\oplus H$ (details
left to the reader). Under this isomorphism, the projection
$H_{\rho'}\to H$ is identified with the map $H_\rho\oplus H\to H$
given by the rule : $(h_1,h_2)\mapsto\pi_H(h_1)\cdot h_2$, where
$\pi_H:H_\rho\to H$ is the projection.
It then follows that the natural map :
$$
Q_\rho\otimes_{H_\rho}H\to Q_{\rho'}\otimes_{H_{\rho'}}H
$$
is an isomorphism. Thus, we may replace $G$ and $\rho$ by $G'$ and
$\rho'$, which allows to assume from start that $\rho$ is surjective.
However, we have natural isomorphisms :
$$
\Ker(Q_\rho\xrightarrow{\ \pi_Q\ }Q)\simeq\Ker\,\rho\simeq
\Ker(H_\rho\xrightarrow{\ \pi_H\ }H).
$$
Moreover, the set underlying $Q$ (resp. $H$) is the set-theoretic
quotient of the set $Q_\rho$ (resp. $H_\rho$) by the translation
action of $\Ker\,\rho$; hence the natural maps
$Q_\rho/\Ker\,\pi_Q\to Q$ and $H_\rho/\Ker\,\pi_H\to H$ are
isomorphisms (lemma \ref{lem_special-p-out}(ii)). We can then compute :
$$
Q_\rho\otimes_{H_\rho}H\simeq
(Q_\rho/\Ker\,\pi_Q)\otimes_{H_\rho/\Ker\,\pi_H}H
\simeq Q\otimes_HH\simeq Q
$$
as stated.
\end{proof}

\begin{proposition}\label{prop_good-charts}
Let $T$ be a locally ringed topos, $\xi$ any $T$-point, and
$\underline M$ a coherent (resp. fine) log structure on $T$.
Suppose that $G$ is a finitely generated abelian group with
a group homomorphism $G\to\underline M{}^\gp_\xi$ such that
the induced map $G\to(\underline M^\sharp)^\gp_\xi$
is surjective. Set
$$
P:=G\times_{\underline M{}^\gp_\xi}\underline M{}_\xi.
$$
Then the induced morphism $P\to\underline M{}_\xi$ extends to
a finite (resp. fine) chart $P_U\to\underline M{}_{|U}$ on some
neighborhood $U$ of\/ $\xi$.
\end{proposition}
\begin{proof} We begin with the following :

\begin{claim}\label{cl_better-chart}
Let $Y$ be any locally ringed topos, $\xi$ a $Y$-point,
and $\alpha:Q_{\!Y}\to\cO_{\!Y}$ the constant pre-log
structure associated to a map of monoids
$\vartheta:Q\to\Gamma(Y,\cO_{\!Y})$, where $Q$ is finitely
generated. Set
$S:=\alpha^{-1}_\xi\cO^\times_{\!Y,\xi}\subset Q_{\!Y,\xi}=Q$.
Then :
\begin{enumerate}
\item
$S$ and $S^{-1}Q$ are finitely generated monoids.
\item
There exists a neighborhood $U$ of $\xi$ such that $\alpha_{|U}$
factors as the composition of the natural map of sheaves of
monoids $j_U:Q_{\!U}\to(S^{-1}Q)_U$, and a (necessarily
unique) pre-log structure $\alpha_S:(S^{-1}Q)_U\to\cO_{\!U}$.
\item
The induced map of log structures
$j_U^{\log}:Q^{\log}_{\!U}\to(S^{-1}Q)^{\log}_U$ is an isomorphism.
\item
$\alpha^{-1}_{S,\xi}(\cO^\times_{\!U,\xi})=(S^{-1}Q)^\times$
is a finitely generated group.
\end{enumerate}
\end{claim}
\begin{pfclaim} (i) follows from lemma \ref{lem_face}(iv).

(ii): Since $\cO^\times_{\!Y,\xi}$ is the filtered colimit of the
groups $\Gamma(U,\cO^\times_{\!U})$, where $U$ ranges over the
neighborhoods of $\xi$, lemma \ref{lem_finite-pres}(ii) and (i)
imply that the induced map $S\to\cO^\times_{\!Y,\xi}$ factors
through $\Gamma(U,\cO^\times_{\!Y})$ for some neighborhood $U$ of
$\xi$. Then the composition of $\vartheta$ and the natural map
$\Gamma(Y,\cO_{\!Y})\to\Gamma(U,\cO_{\!U})$ extends to a unique map
$S^{-1}Q\to\Gamma(U,\cO_{\!U})$, whence the sought pre-log structure
$\alpha_S$ on $U$.

(iii): Let $\underline{N}$ be any log structure on $U$;
it is clear that every morphism of pre-log structures
$Q_U\to\underline{N}$ factors uniquely through $(S^{-1}Q)_U$,
whence the contention.

(iv): Indeed, by construction we have :
$\alpha^{-1}_{S,\xi}(\cO^\times_{\!U,\xi})=S^\gp$.
\end{pfclaim}

Let $Y$ be a neighborhood of $\xi$ such that
$\underline M{}_{|Y}$ admits a finite local chart
$\alpha:Q_{\!Y}\to\cO_{\!Y}$; we lift $\xi$ to some $Y$-point,
$\xi_Y$, and choose a neighborhood $U\in\Ob(T\!/Y)$ of $\xi_Y$,
as provided by claim \ref{cl_better-chart}.
We may then replace $T$ by $U$, $\xi$ by $\xi_Y$ and $\alpha$
by the chart $\alpha_S$ of claim \ref{cl_better-chart}(ii),
which allows to assume from start that
$S:=\alpha^{-1}_\xi(\cO^\times_{\!T,\xi})$ is a finitely
generated group. Moreover, let
$H:=\Ker(S\to\cO^\times_{\!T,\xi})$; clearly
$\alpha_\xi:Q\to\cO_{\!T,\xi}$ factors through the quotient
$Q':=Q/H$, hence we may find a neighborhood $U$ of $\xi$ such
that $\alpha_{|U}$ factors through a (necessarily unique)
map of pre-log structures $\alpha_H:Q'_{\!U}\to\cO_{\!U}$.
Furthermore, if $\underline N$ is any log structure on $U$,
every map of pre-log structures $Q_{\!U}\to\underline N$
factors through $Q'_{\!U}$, so that $\alpha_H$ is a chart
for $\underline M{}_{|U}$.
We may therefore replace $T$ by $U$ and $\alpha$ by $\alpha_H$,
which allows to assume additionally, that $\alpha_\xi$ is
injective on the subgroup $S$. Now, for any finitely generated
subgroup $H\subset\cO^\times_{\!T,\xi}$ with $S\subset H$, we
set $\underline M{}_{\xi,H}:=H\amalg_SQ$; clearly, the
monoids $\underline M{}_{\xi,H}$ are finitely generated, and
moreover :
$$
\underline M{}_\xi=\colim_{S\subset H\subset\cO^\times_{\!T,\xi}}
\underline M{}_{\xi,H}.
$$
Furthermore, we deduce a natural sequence of maps of monoids :
\set\begin{equation}\label{eq_exact-seq-mods}
\{1\}\to\underline M{}_{\xi,H}\xrightarrow{\phi_H}\underline M{}_\xi
\xrightarrow{\psi_H}\cO^\times_{\!T,\xi}/H\to\{1\}.
\end{equation}

\begin{claim}\label{cl_exac-seq-mds}
For every subgroup $H$ as above, the sequence
\eqref{eq_exact-seq-mods} is exact, {\em i.e.}
$\underline M{}_{\xi,H}$ is the kernel of $\psi_H$, and
$\cO^\times_{\!T,\xi}/H$ is the cokernel of $\phi_H$.
\end{claim}
\begin{pfclaim} By lemma \ref{lem_forget-me-not}(iii), the
assertion concerning $\Ker\,\psi_H$ can be verified on the
underlying map of sets; however, lemma \ref{lem_special-p-out}(ii)
says that the set $\underline M{}_\xi$ is the set-theoretic
quotient $(Q\times\cO^\times_{\!T,\xi})/S$, for the natural
translation action of $S$, and a similar
description holds for $\underline{M}_{\xi,H}$, therefore
$\Ker\,\phi_H$ is the set-theoretic quotient $(H\times Q)/S$,
as required.

The assertion concerning $\Coker\,\phi_H$ holds by general
categorical nonsense.
\end{pfclaim}

Let $\eps:\underline M{}_\xi\to\underline M{}^\gp_\xi$ be
the natural map; claim \ref{cl_exac-seq-mds} implies that
the sequence of abelian groups $\eqref{eq_exact-seq-mods}^\gp$
is right exact, and then a little diagram chase shows that :
\set\begin{equation}\label{eq_consequences}
\eps^{-1}(\Img\,\phi_H^\gp)=\underline M{}_{\xi,H}
\qquad
\text{whenever $S\subset H\subset\cO^\times_{\!T,\xi}$}.
\end{equation}
Since $G$ is finitely generated, we may find $H$ as above,
large enough, so that $\underline M{}^\gp_{\xi,H}$ contains
the image of $G$. In view of \eqref{eq_consequences}, we
deduce that the natural map
$$
G\times_{\underline M{}^\gp_{\xi,H}}\underline M{}_{\xi,H}\to P
$$
is an isomorphism, so $P$ is finitely generated, by corollary
\ref{cor_fibres-are-fg}; moreover $P$ is integral whenever
$\underline M{}_\xi$ is. Then, lemma \ref{lem_simple-charts-top}(iv.a)
implies that the natural map $P\to\underline M{}_\xi$ extends
to a morphism of log structures
$\vartheta:P_U^{\log}\to\underline M{}_{|U}$ on some neighborhood
$U$ of $\xi$.
It remains to show that $\vartheta$ restricts to a chart for
$\underline M{}_{|V}$, on some smaller neighborhood $V$ of $\xi$.
To this aim, it suffices to show that the map of stalks
$\vartheta_\xi$ is an isomorphism (lemma
\ref{lem_simple-charts-top}(iv.c)). The latter assertion follows
from lemma \ref{lem_charts-later}.
\end{proof}

Proposition \ref{prop_good-charts} is the basis of several
frequently used tricks that allow to construct ``good'' charts
for a given coherent log structure (and for a morphism of such
structures), or to ``improve'' given charts. We conclude this
section with a selection of these tricks.

\begin{corollary}\label{cor_first-trick}
Let $T$ be a topos, $\xi$ a $T$-point, $\underline M$ a fine
log structure on $T$. Then there exist a neighborhood $U$ of
$\xi$ in $T$, and a chart $P_U\to\underline M{}_{|U}$ such that :
\begin{enumerate}
\alphaenu
\item
$P^\gp$ is a free abelian group of finite rank.
\item
The induced morphism of monoids $P\to\cO_{T,\xi}$ is local.
\end{enumerate}
\end{corollary}
\begin{proof} Choose a group homomorphism
$G:=\Z^{\oplus n}\to\underline M{}^\gp_\xi$ (for some integer
$n\geq 0$), such that the induced map
$G\to(\underline M^\sharp)^\gp_\xi$ is surjective, and set
$P:=G\times_{\underline M{}^\gp_\xi}\underline M{}_\xi$.
By proposition \ref{prop_good-charts}, the induced map
$P\to\underline M{}_\xi$ extends to a chart
$P_U\to\underline M{}_{|V}$, for some neighborhood
$U$ of $\xi$. According to example \ref{ex_regular}(v),
$P^\gp$ is a subgroup of $G$, whence (a). Next, claim
\ref{cl_better-chart} implies that, after replacing $P$
by some localization (which does not alter $P^\gp$), and
$U$ by a smaller neighborhood of $\xi$, we may achieve (b)
as well.
\end{proof}

\begin{corollary}\label{cor_sharpy}
Let $T$ be a topos, $\xi$ a $T$-point, $\underline M$ a coherent
log structure on $T$, and suppose that $\underline M_\xi$ is
integral and saturated. Then we have :
\begin{enumerate}
\item
There exist a neighborhood $U$ of $\xi$ in $T$, and a fine and
saturated chart $P_U\to\underline M{}_{|U}$ which is sharp at the
point $\xi$.
\item
Especially, there exists a neighborhood $U$ of $\xi$ in $T$,
such that $\underline M{}_{|U}$ is a fine and saturated log structure.
\end{enumerate}
\end{corollary}
\begin{proof}(i): By lemma \ref{lem_decomp-sats}, we may find a
decomposition $\underline M{}_\xi=P\times\underline M{}_\xi^\times$,
for a sharp submonoid $P\subset\underline M{}_\xi$. Set $G:=P^\gp$;
we deduce an isomorphism
$G\isom\underline M{}^\gp_\xi/\underline M{}^\times_\xi$, and
clearly $P=G\times_{\underline M{}^\gp_\xi}\underline M{}_\xi$.
By proposition \ref{prop_good-charts}, it follows that the induced
map $P\to\underline M{}_\xi$ extends to a chart
$\beta:P_U\to\underline M{}_{|U}$ on a neighborhood $U$ of
$\xi$. By construction, $\beta$ is sharp at the $T$-point $\xi$;
moreover, since $\underline M{}_\xi$ is saturated, it is easily
seen that the same holds for $P$. Finally, $P$ is finitely
generated, by corollary \ref{cor_fibres-are-fg}.

(ii): This follows immediately from (i) and lemma
\ref{lem_simple-charts-top}(iii).
\end{proof}

\begin{theorem}\label{th_good-charts}
Let $T$ be a locally ringed topos, $\xi$ a $T$-point,
$f:\underline M\to\underline N$ a morphism of coherent
(resp. fine) log structures on $T$. Then :
\begin{enumerate}
\item
There exists a neighborhood $U$ of\/ $\xi$, such that $f_{|U}$
admits a finite (resp. fine) chart.
\item
More precisely, given a finite (resp. fine) chart
$\omega_P:P_T\to\underline M$, we may find a neighborhood
$U$ of\/ $\xi$, a finite (resp. fine) monoid $Q$, and a chart
of $f_{|U}$ of the form
$$
(\omega_{P|U},\omega_Q:Q_U\to\underline N{}_{|U},\theta:P\to Q).
$$
\item
Moreover, if $f$ is a flat (resp. saturated) morphism of fine
log structures and $(\omega_P,\omega_Q,\theta)$ is a fine chart
for $f$, then we may find a neighborhood $U$ of $\xi$, a
localization map $j:Q\to Q'$, and a flat (resp. saturated) and
fine chart for $f_{|U}$ of the form
$(\omega_{P|U},\omega_{Q'},j\circ\theta)$, such that
\begin{enumerate}
\item
$\omega_{Q|U}=\omega_{Q'}\circ j_U$.
\item
The chart $\omega_{Q'}$ is local at the point $\xi$.
\end{enumerate}
\end{enumerate}
\end{theorem}
\begin{proof} Up to replacing $T$ by $T\!/U'_i$ for a covering
$(U'_i\to 1_T~|~i\in I)$ of the final object, we may assume
that we have finite (resp. fine) charts $\omega_P:P_T\to\underline M$
and $Q'_T\to\underline N$ (lemma \ref{lem_simple-charts-top}(iii)),
whence a morphism of pre-log structures :
$$
\omega:P_T\xrightarrow{\ \omega_P\ }
\underline M\xrightarrow{\ f\ }\underline N.
$$
Let $\xi$ be any $T$-point; there follow maps of monoids
$\phi:P\to\underline M{}_\xi\to\underline N_\xi$ and
$\psi:Q'\to\underline N{}_\xi$. Set $G:=(P\oplus Q')^\gp$, and
apply proposition \ref{prop_good-charts} to the group homomorphism
$G\to\underline N{}_\xi^\gp$ induced by $\phi$ and $\psi$; for
$Q:=G\times_{\underline N{}^\gp_\xi}\underline N{}_\xi$, we obtain a
finite (resp. fine) local chart $Q_U\to\underline N{}_{|U}$ on some
neighborhood $U$ of $\xi$. Then, $\phi$ and the natural map
$P\to G$ determine a unique map $P\to Q$, whence a morphism
$\omega':P_U\to Q_U\to\underline N{}_{|U}$ of pre-log structures; by
construction, $\omega'_\xi:P=P_{U,\xi}\to\underline N{}_\xi$ is none
else than $\phi$. By lemma \ref{lem_simple-charts-top}(iv.b), we may
then find a smaller neighborhood $V\to U$ of $\xi$ such that
$\omega'_{|V}=\omega_{|V}$. This proves (i) and (ii).

Next, we suppose that $f$ is flat (resp. saturated) and both
$\underline M$, $\underline N$ are fine, and we wish to show (iii).
In view of claim \ref{cl_better-chart}, we may find -- after
replacing $T$ by a neighborhood of $\xi$ -- a fine chart for
$f$ of the form $(\omega_{P'},\omega_{Q'},\theta')$, such that :
\begin{itemize}
\item
$P'$ and $Q'$ are localizations of $P$ and $Q$, and $\theta'$
is induced by $\theta$;
\item
$\omega_P=\omega_{P'}\circ(j_P)_T$ and
$\omega_Q=\omega_{Q'}\circ(j_Q)_T$, where $j_P:P\to P'$ and
$j_Q:Q\to Q'$ are the localization maps;
\item
the induced maps $Q^{\prime\sharp}\to\underline M{}^\sharp_\xi$
and $P^{\prime\sharp}\to\underline N{}^\sharp_\xi$ are isomorphisms.
\end{itemize}
Now, by proposition \ref{prop_finally-flat} (resp. corollary
\ref{cor_satur-check}), the map $f_\xi$ is flat (resp. saturated),
hence the same holds for the induced map
$\underline M{}^\sharp_\xi\to\underline N{}^\sharp_\xi$, by
corollary \ref{cor_yet-another-flat}(i) (resp. by lemma
\ref{lem_little}(iii)). Then the map
$P^{\prime\sharp}\to Q^{\prime\sharp}$ induced by
$\theta'$ is flat (resp. saturated) as well, so the same
holds for $\theta'$, by corollary \ref{cor_yet-another-flat}(ii)
(resp. again by lemma \ref{lem_little}(iii)). Finally,
$j_P\circ\theta':P\to Q'$ is flat (resp. saturated), by example
\ref{ex_flat-quots}(iii) (resp. by lemma \ref{lem_little}(i)),
and $(\omega_P,\omega_{Q'},j_P\circ\theta')$ is a chart for $f$
with the sought properties.
\end{proof}

\begin{corollary}\label{cor_sharp-chart}
Let $T$ be a topos, $\xi$ a $T$-point, and
$\phi:\underline M\to\underline N$ a saturated morphism of fine
and saturated log structures on $T$. We have :
\begin{enumerate}
\item
There exist a neighborhood $U$ of\/ $\xi$, and a fine and saturate
chart $(\omega_P,\omega_Q,\theta:P\to Q)$ of $\phi_{|U}$, such that
$\omega_P$ and $\omega_Q$ are sharp at the point $\xi$.
\item
More precisely, suppose that
$(\omega_P,\omega_Q,\theta:P\to Q)$ is a fine and saturated chart
for $\phi$, such that $\underline M$ is sharp at the point $\xi$,
and $\omega_Q$ is local at $\xi$. Then there exists a section
$\sigma:Q^\sharp\to Q$ of the projection $Q\to Q^\sharp$, such
that $(\omega_P,\omega_Q\circ\sigma_T,\theta^\sharp)$ is a chart
for $\phi$.
\end{enumerate}
\end{corollary}
\begin{proof} After replacing $T$ by some neighborhood of $\xi$,
we may assume that $\underline M$ admits a chart
$\omega_P:P_T\to\underline M$ which is fine, saturated, and sharp
at the point $\xi$ (corollary \ref{cor_sharpy}(i)). Then, by
theorem \ref{th_good-charts}(iii), we may find a neighborhood
$U$ of $\xi$, and a fine and saturated chart
$(\omega_{P|U},\omega_Q,\theta:P\to Q)$ for $\phi_{|U}$, such
that $\omega_Q$ is local at $\xi$, so that $\theta$ is also a
local morphism. Hence, it suffices to show assertion (ii).

(ii): We notice the following :

\begin{claim}\label{cl_section-sharp}
Let $\theta:P\to Q$ a local and saturated morphism of fine and
saturated monoids, and suppose that $P$ is sharp. Then there
exists a section $\sigma:Q^\sharp\to Q$ of the projection
$Q\to Q^\sharp$, such that $\theta(P)$ lies in the image of $\sigma$.
\end{claim}
\begin{pfclaim} Pick an isomorphism
$\beta:Q\isom Q^\sharp\times Q^\times$ as in lemma
\ref{lem_decomp-sats}, and denote by $\psi:P\to Q^\times$
the composition of $\theta$ with the induced projection
$Q\to Q^\times$. The morphism $\theta^\sharp$ is still local
and saturated (lemma \ref{lem_little}(iii)), hence corollary
\ref{cor_persist-integr}(ii) implies that $\theta^{\sharp\gp}$
extends to an isomorphism $P^\gp\oplus L\isom Q^{\sharp\gp}$,
where $L$ is a free abelian group of finite rank. Thus, we
may extend $\psi^\gp$ to a group homomorphism
$\psi':Q^{\sharp\gp}\to Q^\times$. Define an automorphism
$\alpha$ of $Q^\sharp\times Q^\times$, by the rule :
$(x,g)\mapsto(x,g\cdot\psi'(x)^{-1})$. The restriction
$\sigma:Q^\sharp\to Q$ of $(\alpha\circ\beta)^{-1}$ will do.
\end{pfclaim}

With the notation of claim \ref{cl_section-sharp} it is easily
seen that $\omega_Q\circ\sigma_T$ is still a chart for
$\underline N$, hence $(\omega_P,\omega_Q\circ\sigma_T,\theta^\sharp)$
is a chart for $\phi$ as required.
\end{proof}

\begin{corollary}\label{cor_was-strict-locus-i}
Let $f:(T,\underline M)\to(S,\underline N)$ a morphism of log
topoi with coherent (resp. fine) log structures, and suppose
that $\underline N$ admits a finite (resp. fine) chart
$\omega_Q:Q_S\to\underline N$. Let also $\xi$ be any $T$-point;
we have :
\begin{enumerate}
\item
There exist a neighborhood $U$ of $\xi$, and a finite
(resp. fine) chart $(\omega_{Q|U},\omega_P,\theta:Q\to P)$
for the morphism $f_{|U}$.
\item
Moreover, if $\underline M$ and $\underline N$ are fine and $f$
is log flat (resp. saturated) then, on some neighborhood $U$
of $\xi$, we may also find a chart $(\omega_{Q|U},\omega_P,\theta)$
which is flat (resp. saturated) and fine.
\end{enumerate}
\end{corollary}
\begin{proof} This is an immediate consequence of theorem
\ref{th_good-charts}.
\end{proof}

\subsection{Log schemes}
We specialize now to the case of a scheme $X$. Whereas in
\cite[\S6.4]{Ga-Ra} we considered only pre-log structures on the
Zariski site of a scheme, hereafter we shall treat uniformly the
categories of log structures on the topoi $X^\sim_\et$ and
$X^\sim_\Zar$ (notation of \eqref{subsec_topolo-on-sch}).

\sset\subsubsection{}\label{subsec_special-schs}
Henceforth, we choose $\tau\in\{\et,\Zar\}$ (see
\eqref{subsec_tau-notate}), and whenever we mention a topology
on a scheme $X$, it will be implicitly meant that this is the
topology $X_\tau$ (unless explicitly stated otherwise). Let $X$
be a scheme; a {\em pre-log structure\/} (resp. {\em a log
structure}) on $X$ is a pre-log structure (resp. a log structure)
on the topos $X_\tau^\sim$. The datum of a scheme $X$ and a log
structure on $X$ is called briefly a {\em log scheme}. It is known
that a morphism of schemes $X\to Y$ is the same as a morphism of
locally ringed topoi $X^\sim_\tau\to Y^\sim_\tau$, hence we may
define a morphism of log schemes
$(X,\underline{M})\to(Y,\underline{N})$ as a morphism of log topoi
$(X^\sim_\tau,\underline{M})\to(Y^\sim_\tau,\underline{N})$ (and
likewise for morphisms of schemes with pre-log structures). We
denote by $\bprelog_\tau$ (resp. $\blog_\tau$) the category of
schemes with pre-log structures (resp. of log schemes) on the chosen
topology $\tau$. We denote by :
$$
\bintlog_\tau \qquad \bsatlog_\tau \qquad \bqcohlog_\tau \qquad
\bcohlog_\tau
$$
the full subcategory of the category $\blog_\tau$, consisting of all
log schemes with  integral (resp. integral and saturated, resp.
quasi-coherent, resp. coherent) log structures.

A scheme with a quasi-fine (resp. fine, resp. quasi-fine and saturated,
resp. fine and saturated) log structure is called, briefly, a
{\em quasi-fine log scheme} (resp. a {\em fine log scheme},
resp. a {\em qfs log scheme}, resp. a {\em fs log scheme}),
and we denote by
$$
\bqflog_\tau \qquad \bflog_\tau \qquad \bqfslog_\tau \qquad
\bfslog_\tau
$$
the full subcategory of $\blog_\tau$ consisting of all quasi-fine
(resp. fine, resp. qfs, resp. fs) log schemes on the topology
$\tau$. In case it is clear (or indifferent) which topology we are
dealing with, we will usually omit the subscript $\tau$. There is an
obvious (forgetful) functor :
$$
F:\blog\to\Sch
$$
to the category of schemes, and it is easily seen that this functor
is a fibration. For every scheme $X$, we denote by $\blog_X$ the
fibre category $F^{-1}(X)$ {\em i.e.} the category of all log
structures on $X$ (or $\blog_{X_\tau}$, if we need to specify the
topology $\tau$). The same notation shall be used also for the
various subcategories : so for instance we shall write $\bintlog_X$
for the full subcategory of all integral log structures on $X$.
Moreover, we shall say that the log scheme $(X,\underline M)$ is
{\em locally noetherian\/} if the underlying scheme $X$ is locally
noetherian.

\sset\subsubsection{}\label{subsec_choose-a-top} Most of the
forthcoming assertions hold in both the \'etale and Zariski topoi,
with the same proof. However, it may occasionally happen that the
proof of some assertion concerning $X^\sim_\tau$ (for
$\tau\in\{\et,\Zar\}$), is easier for one choice or the other of
these two topologies; in such cases, it is convenient to be able to
change the underlying topology, to suit the problem at hand. This is
sometimes possible, thanks to the following general considerations.

The morphism of locally ringed topoi $\tilde u$ of
\eqref{subsec_big-site} induces a pair of adjoint functors :
$$
\tilde u{}^*:\blog_\Zar\to\blog_\et \qquad \tilde
u_*:\blog_\et\to\blog_\Zar
$$
as well as analogous adjoint pairs for the corresponding categories
of sheaves of monoids (resp. of pre-log structures) on the two
sites. It follows formally that $\tilde u{}^*$ sends constant log
structures to constant log structures, {\em i.e.} for every scheme
$X$, and every object $M:=(M,\phi)$ of $\Mnd/\Gamma(X,\cO_{\!X})$ we
have a natural isomorphism :
$$
\tilde
u{}^*(X_\Zar,M^{\log}_{X_\Zar})\simeq(X_\et,M^{\log}_{X_\et}).
$$
More generally, lemma \ref{lem_simple-charts-top}(i) shows that
$\tilde u{}^*$ preserves the subcategories of quasi-coherent (resp.
coherent, resp. integral, resp. fine, resp. fine and saturated)
log structures.

\begin{proposition}\label{prop_reduce-to-etale}
{\em (i)}\ \ The functor\/ $\tilde u{}^*$ on log structures is
faithful and conservative.
\begin{enumerate}
\addenu
\item
The functor\/ $\tilde u{}^*$ restricts to a fully faithful functor :
$$
\tilde u{}^*:\bintlog_\Zar\to\bintlog_\et.
$$
\item
Let $(X_\et,\underline M)$ be any log scheme. Then the counit
of adjunction
$\tilde u{}^*\tilde u_*(X_\et,\underline M)\to(X_\et,\underline M)$
is an isomorphism if and only if the same holds for the counit of
adjunction
$\tilde u{}^*\tilde u_*(\underline M^\sharp)\to\underline M^\sharp$.
\end{enumerate}
\end{proposition}
\begin{proof}(i): Let $(X_\Zar,\underline M)$ be any log structure;
set $(X_\et,\underline M{}_\et):= \tilde u{}^*(X,\underline M)$,
and denote by
$\tilde u^\natural_X:\tilde u{}^*\cO_{\!X_\Zar}\to\cO_{\!X_\et}$
the natural map of structure rings. Since $\tilde u$ is a morphism
of locally ringed topoi, we have :
$$
(\tilde u^\natural_X)^{-1}\cO_{\!X_\et}^\times=\tilde
u{}^*(\cO_{\!X_\Zar}^\times).
$$
It follows easily that $\underline M{}_\et$ is the push-out in the
cocartesian diagram :
$$
\xymatrix{ \tilde u{}^*\underline M^\times \ar[r]^-\alpha \ar[d] &
\cO_{\!X_\et}^\times \ar[d] \\ \tilde u{}^*\underline M
\ar[r]^-\beta & \underline M{}_\et. }
$$
However, for every geometric point $\xi$ of $X$, the natural map
$\cO_{\!X,|\xi|}\to\cO_{\!X,\xi}$ is faithfully flat, hence
$\alpha_\xi$ is injective, and therefore also $\beta_\xi$, in light
of lemma \ref{lem_special-p-out}(ii). The faithfulness of $\tilde
u^*$ is an easy consequence. Next, let $f:\underline M\to\underline N$
be a morphism of log structures on $X_\Zar$, set
$(X_\et,\underline N{}_\et):=\tilde u{}^*(X,\underline N)$, and suppose
that $\tilde u{}^*f:\underline M{}_\et\to\underline N{}_\et$ is an
isomorphism; we wish to show that $f$ is an isomorphism. However,
$\beta$ induces an isomorphism of monoids :
\set\begin{equation}\label{eq_induced-by-beta}
\tilde u{}^*(\underline M^\sharp)\isom\underline M{}_\et^\sharp
\end{equation}
and likewise for $\underline N$; it follows already that $f$ induces
an isomorphism
$\underline M^\sharp\isom\underline N^\sharp$.
To conclude, it suffices to invoke lemma \ref{lem_check-iso}.

(ii): Let us suppose that $\underline M$ is integral. According to
\cite[Prop.3.4.1]{Bor}, it suffices to show that the unit of
adjunction $(X,\underline M)\to\tilde u_*(X_\et,\underline M{}_\et)$
is an isomorphism. Now, from the isomorphism \eqref{eq_induced-by-beta}
we deduce a commutative diagram :
$$
\xymatrix{\underline M^\sharp \ar[r] \ar[d] &
(\tilde u_*\underline M{}_\et)^\sharp \ar[d] \\
\tilde u_*\tilde u{}^*(\underline M^\sharp) \ar[r]
& \tilde u_*(\underline M{}_\et^\sharp)}
$$
whose bottom arrow is an isomorphism, and whose left vertical arrow
is an isomorphism as well, by lemma \ref{lem_Hilbert90}(iii) (and
again \cite[Prop.3.4.1]{Bor}). We claim that also the right vertical
arrow is an isomorphism. Indeed, since $\tilde u_*$ is left exact,
the latter arrow is a monomorphism, hence it suffices to show it is
an epimorphism; however, since $\underline M{}_\et$ is an integral log
structure (lemma \ref{lem_simple-charts-top}(i)), it is easily seen
that the projection $\underline M{}_\et\to\underline M{}_\et^\sharp$ is
a $\cO_{\!X_\et}^\times$-torsor (in the topos
$X^\sim_\et/\underline M{}_\et^\sharp$). Then the contention follows
from the exact sequence of pointed sheaves \eqref{eq_non-ab-coh},
and the vanishing result of lemma \ref{lem_Hilbert90}(iv).

Summing up, we conclude that the top horizontal arrow in the above
diagram is an isomorphism, so the assertion follows from lemma
\ref{lem_check-iso}.

(iii): Set $(X_\et,(\tilde u_*\underline M)_\et):=
\tilde u{}^*\tilde u_*(X_\et,\underline M)$. To begin with,
lemma \ref{lem_Hilbert90}(iv) easily implies that the natural
morphism
$(\tilde u_*\underline M)^\sharp\to\tilde u_*(\underline M^\sharp)$
is an isomorphism; together with the general isomorphism
\eqref{eq_general-iso}, this yields a short exact sequence
of $X_\et$-monoids :
$$
0\to\cO^\times_{\!X_\et}\to(\tilde u_*\underline M)_\et\to
\tilde u{}^*\tilde u_*(\underline M^\sharp)\to 0
$$
which easily implies the assertion : the details shall be
left to the reader.
\end{proof}

We shall prove later on some more results in the same vein (see
corollary \ref{cor_coherence-goes-down}).

\sset\subsubsection{}
Arguing as in \eqref{subsec_morph-logtop}, we see easily
that all finite limits are representable in the category of
log schemes. The rule $X\mapsto(X,\cO_{\!X}^\times)$ defines
a fully faithful inclusion of the category of schemes into the
category of log schemes. Hence, we shall regard a scheme as a
log scheme with trivial log structure. Especially, if
$(X,\underline M)$ is any log scheme, and $Y\to X$ is any
morphism of schemes, we shall often use the notation :
\set\begin{equation}\label{eq_prod-sch-log-sch}
Y\times_X(X,\underline M):=
(Y,\cO_Y^\times)\times_{(X,\cO_{\!X}^\times)}(X,\underline M).
\end{equation}
Especially, if $\xi$ is any $\tau$-point of $X$, the
{\em localization} of $(X,\underline M)$ at $\xi$ is the
log scheme
$$
(X(\xi),\underline M(\xi)):=X(\xi)\times_X(X,\underline M)
$$
(see definition \ref{def_strict-loc}(ii,iii)).
If $\tau=\et$, this operation is also called the
{\em strict henselization} of $X$ at $\xi$.

\begin{definition}\label{def_trivial-locus}
(i)\ \
For every integer $n\in\N$, we have the subset :
$$
(X,\underline M)_n:=\{x\in X~|~\dim\underline M_\xi\leq n\ 
\text{for every $\tau$-point $\xi\to X$ localized at $x$}\}.
$$
Especially $(X,\underline M)_0$ consists of all points $x\in X$
such that $\underline M_\xi=\cO^\times_{\!X,\xi}$ for every
$\tau$-point $\xi$ of $X$ localized at $x$; this subset is called
the {\em trivial locus\/} of $(X,\underline M)$, and is also denoted
$(X,\underline M)_\tr$.

(ii)\ \
If $f:(X,\underline M)\to(Y,\underline N)$
is a morphism of log schemes, we denote by
$\mathrm{Str}(f)\subset X$ the {\em strict locus\/} of $f$,
which is the subset consisting of all points $x\in X$ such that
$f$ is strict at every $\tau$-point localized at $x$ (see
definition \ref{def_amorph-log}(ii)).
\end{definition}

\begin{remark}\label{rem_trivial-is-strict}
Let $f:(X,\underline M)\to(Y,\underline N)$ be any morphism
of log schemes. 

(i)\ \
Since $\log f$ induces local morphisms on stalks, it is
easily seen that $f$ restricts to a map
$$
f_\tr:(X,\underline M)_\tr\to(Y,\underline N)_\tr.
$$

(ii)\ \
Especially, we have $(X,\underline M)_\tr\subset\mathrm{Str}(f)$.
\end{remark}

\sset\subsubsection{}\label{subsec_annoying-general}
Let $\underline X:=(X_i~|~i\in I)$ be a cofiltered family
of quasi-separated schemes, with affine transition morphisms
$f_\phi:X_j\to X_i$, for every morphism $\phi:j\to i$ in $I$.
Denote by $X$ the limit of\/ $\underline X$, and for each
$i\in I$ let $\pi_i:X\to X_i$ be the natural projection.

\begin{lemma}\label{lem_annoying-gen}
In the situation of \eqref{subsec_annoying-general}, let
$\cX:=((X_i,\underline M_i)~|~i\in I)$ be a cofiltered
system of log schemes, with transition morphisms
$(f_\phi,\log f_\phi):(X_i,\underline M_i)\to(X_j,\underline M_j)$
for every morphism $\phi:i\to j$ in $I$. We have :
\begin{enumerate}
\item
The limit of the system $\cX$ exists in the category\/ $\blog$.
\item
Let $(X,\underline M)$ denote the limit of the system $\cX$.
If $X_i$ is quasi-compact for every $i\in I$, then the natural
map
$$
\colim_{i\in I}\Gamma(X_i,\underline N_i)\to\Gamma(X,\underline N)
$$
is an isomorphism.
\end{enumerate}
\end{lemma}
\begin{proof}(i): Let $X$ be the limit of the system of schemes
$\underline X$, and endow $X$ with the sheaf of monoids
$\underline M:=\colim_{i\in I}\pi_i^*\underline M_i$, where
$\pi^*$ is the pull-back functor for sheaves of monoids
(see \eqref{subsec_pts-of-top}), and the transition maps
$\pi^*_j\underline M_j\to\pi_i^*\underline M_i$ are induced
by the morphisms
$\log f_\phi:f_\phi^*\underline M_j\to\underline M_i$, for every
$\phi:i\to j$ in $I$. Then the structure maps of the log
structures $\underline M_i$ induce a well defined morphism
of $X$-monoids $\underline M\to\cO_{\!X}$, and we claim
that the resulting scheme with pre-log structure
$(X,\underline M)$ is actually a log scheme. Indeed, the
assertion can be checked on the stalks, and notice that,
for every $\tau$-point $\xi$ of $X$ we have a natural
identification
$$
\cO_{\!X,\xi}\isom\colim_{i\in I}\cO_{\!X_i,\pi_i(\xi_i)}.
$$
(This is clear for $\tau=\Zar$, and for $\tau=\et$ one uses
\cite[Ch.IV, Prop.18.8.18(ii)]{EGA4}); it then suffices to invoke
lemma \ref{lem_concern}(i). Lastly, it is easily seen that
$(X,\underline M)$ is a limit of the system $\cX$ : the details
shall be left to the reader.

(ii): In view of the explicit construction in (i), the
assertion follows immediately from proposition
\ref{prop_dir-im-and-colim}.
\end{proof}

\begin{example}\label{ex_norm-cross} Let $X$ be a scheme.
For the following example we choose to work with the \'etale
topology $X_\et$ on $X$.
A {\em divisor\/} on $X$ is a closed subscheme $D\subset X$
which is regularly embedded in $X$ and of codimension $1$
(\cite[Ch.IV, D{\'e}f.19.1.3, \S21.2.12]{EGA4}). Suppose moreover
that $X$ is noetherian, let $D\subset X$ be a divisor, and
denote by $(D_i~|~i\in I)$ the irreducible reduced components
of $D$. We say that $D$ is a {\em strict normal crossings\/}
divisor, if :
\begin{itemize}
\item
$\cO_{\!X,x}$ is a regular ring, for every $x\in D$.
\item
$D$ is a reduced subscheme.
\item
For every subset $J\subset I$, the (scheme theoretic)
intersection $\bigcap_{j\in J}D_j$ is regular of pure
codimension $\sharp J$ in $X$.
\end{itemize}
A closed subscheme $D$ of a noetherian scheme $X$ is a
{\em normal crossings divisor\/} if, for every $x\in X$
there exists an {\'e}tale neighborhood $f:U\to X$ of $x$
such that $f^{-1}D$ is a strict normal crossings divisor
in $U$.

Suppose that $D$ is a normal crossings divisor of a noetherian scheme
$X$, and let $j:U:=X\!\setminus\! D\to X$ be the natural open
immersion. We claim that the log structure $j_*\cO^\times_{\!U}$ is
fine (this is the direct image of the trivial log structure on
$U_\et$ : see example \ref{ex_push-trivial}(ii)). To see this, let
$\xi$ be any geometric point of $X$ localized at a point of $D$; up
to replacing $X$ by an {\'e}tale neighborhood of $\xi$, we may
assume that $D$ is a strict normal crossings divisor; we can assume
as well that $X$ is affine and small enough, so that the irreducible
components $(D_\lambda~|~\lambda\in\Lambda)$ are of the form
$V(I_\lambda)$, where $I_\lambda\subset A:=\Gamma(X,\cO_{\!X})$ is a
principal divisor, say generated by an element $x_\lambda\in A$, for
every $\lambda\in\Lambda$. We claim that $j_*\cO^\times_{\!U}$ is
the constant log structure associated to the pre-log structure :
$$
\alpha:\N^{(\Lambda)}_X\to\cO_{\!X}\quad : \quad e_\lambda\mapsto x_\lambda
$$
where $(e_\lambda~|~\lambda\in\Lambda)$ is the standard basis
of $\N^{(\Lambda)}$. Indeed, let $S\subset\Lambda$ be the largest
subset such that the image of $\xi$ lies in
$D_S:=\bigcap_{\lambda\in S}D_\lambda$, we have
$x_\lambda\in\cO^\times_{\!X,x}$ for all $\lambda\notin S$,
so that the push-out of the induced diagram of stalks
$\cO^\times_{\!X,\xi}\leftarrow\alpha^{-1}\cO^\times_{\!X,\xi}
\to\N^{(\Lambda)}$ is the same as the push-out $P_S$ of the analogous diagram
$\cO^\times_{\!X,\xi}\leftarrow\alpha_S^{-1}\cO^\times_{\!X,\xi}\to\N^{(S)}$,
where $\alpha_S:\N^{(S)}_X\to\cO_{\!X}$ is the restriction of $\alpha$.
Suppose that $a\in\cO_{\!X,\xi}$ and $a$ is invertible on
$X(\xi)\!\setminus\! D_S$; the minimal associated primes of $A/(a)$
are all of height one, and they must therefore be found among the
prime ideals $Ax_\lambda$, with $\lambda\in S$. It follows easily
that $a$ is of the form
$u\cdot\prod_{\lambda\in S}x_\lambda^{k_\lambda}$ for certain
$k_\lambda\in\N$ and $u\in\cO^\times_{\!X,\xi}$. Therefore, the
natural map $\beta_\xi:P_S\to(j_*\cO^\times_{\!U})_\xi$ is surjective.
Moreover, the family $(x_\lambda~|~\lambda\in S)$ is a regular
system of parameters of $\cO_{\!X,\xi}$
(\cite[Ch.0, Prop.17.1.7]{EGAIV}), therefore the natural
map $\Sym^n_{\kappa(\xi)}(\fm_\xi/\fm_\xi^2)\to\fm^n_\xi/\fm^{n+1}_\xi$
is an isomorphism for every $n\in\N$ (here $\fm_\xi\subset\cO_{\!X,\xi}$
is the maximal ideal); it follows easily that $\beta_\xi$ is also
injective.
\end{example}

\begin{example} Suppose that $X$ is a regular noetherian scheme,
and $D\subset X$ a divisor on $X$; let $U:=X\!\setminus\! D$.
If $D$ is not a normal crossings divisor, the log structure
$\underline{M}:=j_*\cO_{\!U}^\times$ on $X_\et$ is not necessarily
fine. For a counterexample, let $K$ be an algebraically closed
field, $C\subset\A^2_K$ a nodal cubic; take $D\subset X:=\A^3_K$
to be the (reduced, affine) cone over the cubic $C$, with vertex
$x_0\in\A^3$, and pick a geometric point $\xi$ localized at $x_0$.
It is easily seen that, away from the vertex, $D$ is a normal
crossings divisor, hence $\underline{M}_{|X\setminus\{x_0\}}$
is a fine log structure on $X\!\setminus\!\{x_0\}$, by example
\ref{ex_norm-cross}.
More precisely, let $y_0\in C$ be the unique singular point,
$L\subset D$ the line spanned by $x_0$ and $y_0$, and $\eta$
a geometric point localized at the generic point of $L$.
By inspecting the argument in example \ref{ex_norm-cross}, we
find that :
$$
\underline{M}_\eta\simeq\N^{\oplus 2}\oplus\cO^\times_{\!X,\eta}
$$
(indeed, an isomorphism is obtained by choosing $a,b\in\cO_{C,y_0}$
such that $V(a)$ and $V(b)$ are the two branches of the cubic $C$
in an {\'e}tale neighborhood of $y_0$). On the other hand, let
$\fp\subset A:=K[T_1,T_2,T_3]$ be the prime ideal corresponding to
$x_0$, and $I\subset A_\fp$ the ideal defining the closed subscheme
$X(x_0)\cap D$ in $X(x_0)$; we claim that $I\cdot\cO_{\!X,\xi}$ is
still a prime ideal, necessarily of height one. Indeed, let
$B:=A_\fp/I$, and denote by $A^\wedge$ (resp. $B^\wedge$) the
$\fp$-adic completion of $A_\fp$ (resp. of $B$); then $B^\wedge$ is
also the completion of the reduced ring $\cO_{\!X,\xi}/I$, hence it
suffices to show that $\Spec\,B^\wedge$ is irreducible. However, we
may assume that $C\subset\Spec\,K[T_1,T_2]$ is the affine cubic
defined by the ideal $J\subset K[T_1,T_2]$ generated by
$T_1^3-T^2_2+T_1T_2$. Then $I$ is generated by the element
$f:=T_1^3-T_2^2T_3+T_1T_2T_3$; also, $\fp=(T_1,T_2,T_3)$, so that
$A^\wedge\simeq K[[T_1,T_2,T_3]]$ and $B^\wedge\simeq A^\wedge/(f)$.
Suppose $\Spec\,B^\wedge$ is not irreducible. This means that
$V(f)\subset\Spec\,A^\wedge$ is a union $V(f)=Z_1\cup\cdots\cup Z_n$
of $n\geq 2$ irreducible components $Z_i$. Since $A^\wedge$ is a
local regular ring, each such irreducible component $Z_i$ is a
divisor, defined by some principal prime ideal $\fq_i$ in
$A^\wedge$. Let $a_i$ be a generator for $\fq_i$; then $f$ admits
factorizations of the form $f=a_ib_i$ for some non-invertible
$b_i\in A^\wedge$. Fix some $i$, and set $a:=a_i$, $b:=b_i$; since
$f$ is homogeneous of degree $3$, we must have
$a\in\fp^k\!\setminus\!\fp^{k+1}$ for either $k=1$ or $k=2$ ($k\neq
0$ since $a$ is not a unit, and $k\neq 3$, since $b$ is not a unit);
then $b\in\fp^{3-k}\!\setminus\!\fp^{4-k}$. Write $a=a'+a''$ and
$b=b'+b''$, where $a''\in\fp^{k+1}$, $b''\in\fp^{4-k}$ and $a'$
(resp. $b'$) is homogeneous of degree $k$ (resp. $3-k$). Then
$f=ab=a'b'+c$, where $c\in\fp^4$ and $a'b'$ is homogeneous of degree
$3$. This means that $f=a'b'$ is a factorization of $f$ in $A$.
However, $f$ is irreducible in $A$, a contradiction. (Instead of
this elementary argument, one can appeal to \cite[Th.43.20]{Na},
which runs as follows. If $R$ is a local domain, then there is a
natural bijection between the set of minimal prime ideals of the
henselization $R^\he$ of $R$ and the set of maximal ideals of the
normalization of $R$ in its ring of fractions. In our case, the
normalization $D^\nu$ of $D$ is the cone over the normalization of
$C$, hence the only point of $D^\nu$ lying over $x_0$ is the vertex
of $D^\nu$.)

It follows that any choice of a generator of $I$ yields an
isomorphism :
$$
\underline{M}_\xi\simeq\N\oplus\cO^\times_{\!X,\xi}.
$$
Suppose now that -- in an {\'e}tale neighborhood $V$ of $\xi$ --
the log structure $\underline{M}$ is associated to a pre-log
structure $\alpha:P_V\to\cO_{\!V}$, for some monoid $P$; hence
$\underline{M}_{|V}$ is the push-out of the diagram
$\cO^\times_{\!V}\leftarrow\alpha^{-1}(\cO^\times_{\!V})\to P_V$,
whence isomorphisms :
$$
P/\alpha^{-1}_\xi(\cO^\times_{\!X,\xi})\simeq
\underline{M}_\xi/\cO^\times_{\!X,\xi}\simeq\N
\qquad
P/\alpha^{-1}_\eta(\cO^\times_{\!X,\eta})\simeq
\underline{M}_\eta/\cO^\times_{\!X,\eta}\simeq\N^{\oplus 2}.
$$
But clearly $\alpha^{-1}_\xi(\cO^\times_{\!X,\xi})
\subset\alpha^{-1}_\eta(\cO^\times_{\!X,\eta})$,
so we would have a surjection of monoids $\N\to\N^{\oplus 2}$,
which is absurd.

On the other hand, we remark that the log structure
$j_*\cO^\times_{\!U}$ {\em on the Zariski site\/} $X_\Zar$
is fine : indeed, one has a global chart
$\N_{X_\Zar}\to j_*\cO^\times_{\!U}$, provided by
the equation defining the divisor $D$.
\end{example}

\sset\subsubsection{}\label{subsec_constant-log}
Let $R$ be a ring, $M$ a monoid, and set $S:=\Spec\,R[M]$.
The unit of adjunction $\eps_M:M\to R[M]$ can be regarded
as an object $(M,\eps_M)$ of $\Mnd/\Gamma(S,\cO_{\!S})$,
whence a constant log structure $M^{\log}_{\!S}$ on $S$
(see \eqref{subsec_Konstant}). The rule
$$
M\mapsto\Spec(R,M):=(S,M^{\log}_{\!S})
$$
is clearly functorial in $M$. Namely, to any morphism
$\lambda:M\to N$ of monoids, we attach the morphism of
log schemes
$$
\Spec(R,\lambda):=(\Spec\,R[\lambda],\lambda^{\log}_{\Spec\,R[N]}):
\Spec(R,N)\to\Spec(R,M).
$$
Likewise, if $P$ is a pointed monoid, $\Spec\,R\La P\Ra$
is a closed subscheme of $\Spec\,R[P]$, and we may define
$$
\Spec\La R,P\Ra:=
\Spec(R,P)\times_{\Spec\,R[P]}\Spec\,R\La P\Ra.
$$
Lastly, if $M$ is a non-pointed monoid, notice the
natural isomorphism of log schemes
$$
\Spec\La R,M_\circ\Ra\isom\Spec(R,M)_\circ.
$$
(Notation of remark \ref{rem_strict-locus}(ii) : the
details shall be left to the reader.)

\begin{lemma}\label{lem_localize-const}
With the notation of \eqref{subsec_constant-log}, let $a\in M$
be any element, and set $M_a:=S_a^{-1}M$, where
$S_a:=\{a^n~|~n\in\N\}$. Then $U_a:=\Spec\,R[M_a]$ is an open
subscheme of $S$, and the induced morphism of log schemes :
$$
\Spec(R,M_a)\to\Spec(R,M)\times_SU_a
$$
is an isomorphism.
\end{lemma}
\begin{proof} Let $\beta_S:M_S\to\cO_{\!S}$ and
$\beta_{U_a}:(M_a)_{U_a}\to\cO_{\!U_a}$ be the natural charts,
and denote by $\phi:M\to M_a$ the localization map.
For every $\tau$-point $\xi$ of $U_a$, we have the
identity :
$$
\beta_{S,\xi}=\beta_{U_a,\xi}\circ\phi:M\to\cO_{\!S,\xi}.
$$
Let $Q:=\beta_{U_a,\xi}^{-1}\cO^\times_{\!U_a,\xi}$.
The assertion is a straightforward consequence of
the following :

\begin{claim} The induced commutative diagram of monoids :
$$
\xymatrix{ \phi^{-1}Q \ar[r] \ar[d] & M \ar[d] \\
           Q \ar[r] & M_a
}$$
is cocartesian.
\end{claim}
\begin{pfclaim}[] Let $b\in M$, and suppose that
$\beta_{U_a,\xi}(a^{-1}b)\in\cO^\times_{\!U_a,\xi}$.
Since $\beta_{U_a,\xi}(a)\in\cO^\times_{\!U_a,\xi}$,
we deduce that the same holds for $\beta_{U_a,\xi}(b)$,
{\em i.e.} $b\in\phi^{-1}Q$. Let $Q':=\phi(\phi^{-1}Q)$;
we conclude that $Q=S_a^{-1}Q'$, the submonoid of $M_a$
generated by $Q'$ and $a^{-1}$. The claim follows easily.
\end{pfclaim}
\end{proof}

\sset\subsubsection{}\label{subsec_from-X-to-cst}
In the same vein, let $X$ be a $R$-scheme, and $(M,\phi)$ any
object of $\Mnd/\Gamma(X,\cO_{\!X})$. The map $\phi$ induces,
via the adjunction of \eqref{subsec_mon-to-algs}, a homomorphism
of $R$-algebras $R[M]\to\Gamma(X,\cO_{\!X})$, whence a map of
schemes $f:X\to S:=\Spec\,R[M]$, inducing a morphism
\set\begin{equation}\label{eq_spec-log}
X\times_S\Spec(R,M)\to(X,(M,\phi)_X^{\log})
\end{equation}
of log schemes.

\begin{lemma}\label{lem_time}
In the situation of \eqref{subsec_from-X-to-cst}, we have :
\begin{enumerate}
\item
The map \eqref{eq_spec-log} is an isomorphism.
\item
The log scheme $\Spec(R,M)$ represents the functor
$$
\blog\to\Set
\quad :\quad
(Y,\underline N)\mapsto
\Hom_{\Z\Alg}(R,\Gamma(Y,N))\times\Hom_\Mnd(P,\Gamma(Y,\underline N)).
$$
\end{enumerate}
\end{lemma}
\begin{proof}(i): The log structure $f^*(M_{\!S}^{\log})$ of
$\Spec(R,M)\times_SX$ represents the functor :
$$
F:\blog_X\to\Set\quad:\quad
\underline{N}\mapsto
\Hom_{\Mnd/\Gamma(S,\cO_{\!S})}((M,\eps_M),\Gamma(S,f_*\underline{N})).
$$
However, if $\underline{N}$ is any log structure on $X$, the pre-log
structure $(f_*\underline{N})^\prelog$ is the same as
$f_*(\underline{N}^\prelog)$ (see \cite[(6.4.8)]{Ga-Ra}). From the
explicit construction of direct images for pre-log structures, and
since the global sections functor is left exact (because it is a
right adjoint), we deduce a cartesian diagram of monoids :
$$
\xymatrix{
\Gamma(S,f_*\underline{N}) \ar[r] \ar[d] & \Gamma(S,\cO_{\!S}) \ar[d] \\
\Gamma(X,\underline{N}) \ar[r] & \Gamma(X,\cO_{\!X}).
}$$
It follows easily that $F$ is naturally isomorphic to
the functor given by the rule :
$$
\underline{N}\mapsto
\Hom_{\Mnd/\Gamma(X,\cO_{\!X})}((M,\phi),\Gamma(X,\underline{N})).
$$
The latter is of course the functor represented by
$(M,\phi)^{\log}_{\!X}$.

(ii) can now be deduced formally from (i), or proved directly
by inspecting the definitions.
\end{proof}

\sset\subsubsection{}
From \eqref{subsec_constant-log} it is also clear that the
rule $(R,M)\mapsto\Spec(R,M)$ defines a functor
$$
\Z\Alg^o\times\Mnd^o\to\blog
$$
which commutes with fibre products; namely, say that
$$
(R',M')\leftarrow(R,M)\to(R'',M'')
$$
are two morphisms of $\Z\Alg\times\Mnd$; then there is
a natural isomorphism of log schemes :
$$
\Spec(R'\otimes_RR'',M'\otimes_MM'')\isom
\Spec(R',M')\times_{\Spec(R,M)}\Spec(R'',M'').
$$
For the proof, one compares the universal properties characterizing
these log schemes, using lemma \ref{lem_time}(ii) : details left
to the reader.

\sset\subsubsection{}\label{subsec_rank-functions}
Let $X$ be a scheme, $\underline M$ a sheaf of monoids on
$X_\tau$. We say that $\underline M$ is {\em locally constant\/}
if there exist a covering family
$(U_\lambda\to X~|~\lambda\in\Lambda)$ and for every
$\lambda\in\Lambda$, a monoid $P_\lambda$ and an isomorphism
of sheaves of monoids
$\underline M{}_{|U_\lambda}\simeq(P_\lambda)_{U_\lambda}$.
We say that $\underline M$ is {\em constructible\/} if,
for every affine open subset $U\subset X$ we can find finitely
many constructible subsets $Z_1,\dots,Z_n\subset U$ such that :
\begin{itemize}
\item
$U=Z_1\cup\cdots\cup Z_n$
\item
$\underline M{}_{|Z_i}$ is a locally constant sheaf of monoids,
for every $i=1,\dots,n$.
\end{itemize}
If moreover, $\underline M{}_\xi$ is a finitely generated monoid
for every $\tau$-point $\xi$ of $X$, we say that $\underline M$
is {\em of finite type}.
If $\phi:\underline M\to\underline N$ is a morphism of constructible
sheaves of monoids on $X_\tau$, then it is easily seen that
$\Ker\,\phi$, $\Coker\,\phi$ and $\Img\,\phi$ are also constructible.
Moreover, if $\underline M$ and $\underline N$ are of finite
type, then the same holds for $\Coker\,\phi$ and $\Img\,\phi$.

Suppose that $\underline M$ is a sheaf of monoids on $X_\tau$,
and for every $x\in X$ choose a $\tau$-point $\bar x$ localized
at $x$; the {\em rank\/} of $\underline M$ is the function 
$$
\rk_{\underline M}:X\to\N\cup\{\infty\}
\qquad
x\mapsto\dim_\Q\underline M{}^\gp_{\bar x}\otimes_\Z\Q.
$$
It is clear from the definitions that the rank function of a
constructible sheaf of monoid of finite type is constructible on $X$.

\begin{lemma}\label{lem_simple-charts}
Let $X$ be a scheme, $\phi:\underline Q\to\underline Q'$ a morphism
of coherent log structures on $X_\tau$. Then :
\begin{enumerate}
\item
The sheaves\/ $\underline Q^\sharp$ and $\Coker\,\phi$ are
constructible of finite type.
\item
$(X,\underline Q)_n$ is an open subset of $X$, for every integer
$n\geq 0$. (See definition {\em\ref{def_trivial-locus}(i)}.)
\item
The rank functions of $\underline Q^\sharp$ and $\Coker\,\phi$
are (constructible and) upper semicontinuous.
\end{enumerate}
\end{lemma}
\begin{proof} Suppose that $\underline{Q}$ admits a finite chart
$\alpha:M_X\to\underline{Q}$. Pick a finite system of generators
$\Sigma\subset M$, and for every $S\subset\Sigma$, set :
$$
Z_S:=\bigcap_{s\in S}D(s)\cap\bigcap_{t\in\Sigma\setminus S}V(t)
$$
where, as usual $D(s)$ (resp. $V(s)$) is the open (resp. closed)
subset of the points $x\in X$ such that the image $s(x)\in\kappa(x)$
of $s$ is invertible (resp. vanishes). Clearly each $Z_S$ is a
constructible subset of $X$, and their union equals $X$. Moreover,
for every $S\subset\Sigma$, and every $\tau$-point $\xi$ supported
on $Z_S$, the submonoid
$N_\xi:=\alpha^{-1}_\xi(\cO^\times_{\!X,\xi})\subset M$ is a face of
$M$ (lemma \ref{lem_face}(i)), hence it is the submonoid $\langle
S\rangle$ generated by $\Sigma\cap N_\xi=S$ (lemma
\ref{lem_face}(ii)). It follows easily that
$\underline Q{}^\sharp_{|Z_S}\simeq(M/\langle S\rangle)_{|Z_S}$.

More generally, suppose that $\underline{Q}$ is coherent. We may
assume that $X$ is quasi-compact. Then, by the foregoing, we may
find a finite set $\Lambda$ and a covering family
$(f_\lambda:U_\lambda\to X~|~\lambda\in\Lambda)$ of $X$, such that
$\underline Q{}^\sharp_{|U_\lambda}$ is a constructible sheaf of
monoids on $U_\lambda$.
Since $f_\lambda$ is finitely presented, it maps constructible
subsets to constructible subsets (\cite[Ch.IV, Th.1.8.4]{EGA4});
it follows easily that the restriction of $\underline Q^\sharp$
to $f_\lambda(U_\lambda)$ is constructible, therefore
$\underline Q^\sharp$ is constructible. Next, notice that
$\Coker\,\phi=\Coker\,\phi^\sharp$; by the foregoing,
$\underline Q^{\prime\sharp}$ is constructible as well, so
the same holds for $\Coker\,\phi$.

Next, Let $x\in X$ be any point, and $\xi$ a $\tau$-point
of $X$ localized at $x$; by theorem \ref{th_good-charts}(i) we
may find a neighborhood $f:U\to X$ of $\xi$ in $X_\tau$ and
a finite chart $(\omega_P,\omega_{P'},\theta)$ for $\phi_{|U}$.
By claim \ref{cl_better-chart}, we may assume, after replacing
$U$ by a smaller neighborhood of $\xi$, that $\omega_P$ and
$\omega_{P'}$ are local at the point $\xi$, in which case
$\underline Q{}^\sharp_\xi=P^\sharp$ and
$\Coker\,\phi_\xi=\Coker\,\theta$. Let $r:X\to\N$ (resp.
$r':U\to\N$) denote the rank function of $\underline Q^\sharp$
(resp. of $\underline Q{}^\sharp_{|U}$); then it is clear that
$r'=r\circ f$. On the other hand, for every $y\in U$ and every
$\tau$-point $\eta$ of $U$ localized at $y$, the stalk
$\underline Q{}^\sharp_\eta$ is a quotient of $P^\sharp$, hence
$r'(y)\leq r(x)$ and
$\dim\underline Q{}_\eta\leq\dim\underline Q{}_\xi$. Since $f$
is an open mapping, this shows (ii), and also that the rank
function of $\underline Q^\sharp$ is upper semicontinuous.
Likewise, let $s$ (resp. $s'$) denote the rank function of
$\Coker\,\phi$ (resp. $\Coker\,\phi_{|U}$); then $s'=s\circ f$,
and $\Coker\,\phi_\eta$ is a quotient of $\Coker\,\theta$ for
every $\tau$-point $\eta$ of $U$; the latter implies that
$s'(y)\leq s(x)$ for every $y\in U$, which shows that $s$ is
upper semicontinuous.
\end{proof}

\begin{corollary}\label{cor_coherence-goes-down}
Let $X$ be a scheme, $\underline M$ a log structure on $X_\Zar$,
and set
$$
(X_\et,\underline M{}_\et):=\tilde u{}^*(X_\Zar,\underline M).
$$
Then $\underline M$ is integral (resp. coherent, resp. fine,
resp. fine and saturated) if and only if the same holds for
$\underline M{}_\et$.
\end{corollary}
\begin{proof} It has already been remarked that
$(X_\et,\underline M{}_\et)$ is coherent (resp. fine, resp. fine and
saturated) whenever the same holds for $(X_\Zar,\underline M)$;
furthermore, the proof of proposition \ref{prop_reduce-to-etale}(i)
shows that the natural map on stalks
$\underline M{}_{|\xi|}\to\underline M{}_{\et,\xi}$ is injective for
every geometric point $\xi$ of $X$, therefore $\underline M$ is
integral whenever the same holds for $\underline M{}_\et$.

Next, we suppose that $\underline M{}_\et$ is coherent, and we
wish to show that $\underline M$ is coherent.

Let $x\in X$ be any point, $\xi$ a geometric point localized at $x$.
By assumption there exists a finitely generated monoid $P'$, with a
morphism $\alpha:P'\to\underline M{}_{\et,\xi}$ inducing an
isomorphism :
$$
P'\otimes_{\beta^{-1}\underline M{}^\times_{\et,\xi}}
\underline M{}^\times_{\et,\xi}\to\underline M{}_{\et,\xi}\simeq
\underline M{}_x\otimes_{\underline M{}^\times_x}\cO^\times_{\!X,\xi}.
$$
It follows easily that we may find a finitely generated submonoid
$Q\subset\underline M{}_x$, such that the image of $\alpha$ lies in
$(Q\cdot\underline M{}_x^\times)\otimes_{\underline M{}^\times_x}
\cO^\times_{\!X,\xi}$, and therefore the natural map :
$$
(Q\cdot\underline M{}_x^\times)\otimes_{\underline M{}^\times_x}
\cO^\times_{\!X,\xi}\to\underline M{}_{\et,\xi}
$$
is surjective. Then lemma \ref{lem_special-p-out}(ii) implies that
$Q\cdot\underline M{}_x^\times=\underline M{}_x$, in other words,
the induced map
$Q\to\underline M{}_x^\sharp= \underline M{}_{\et,\xi}^\sharp$ is
surjective. Set
$$
P:=Q^\gp\times_{\underline M{}^\gp_{\et,\xi}}\underline M{}_{\et,\xi}.
$$
In this situation, proposition \ref{prop_good-charts} tells us that
the induced map $\beta_\xi:P\to\underline M{}_{\et,\xi}$ extends to
an isomorphism of log structures
$\beta^{\log}:P^{\log}_{\!U_\et}\to\underline M{}_{\et|U_\et}$ on
some \'etale neighborhood $U\to X$ of $\xi$. On the other hand, it
is easily seen that the diagram of monoids :
$$
\xymatrix{
\underline M{}_x \ar[r] \ar[d] & \underline M{}^\gp_x \ar[d] \\
\underline M{}_{\et,\xi} \ar[r] & \underline M{}^\gp_{\et,\xi} }
$$
is cartesian, therefore $\beta_\xi$ factors uniquely through a
morphism $\beta'_x:P\to\underline M{}_x$. The latter extends to a
morphism of log structures
$\beta^{\prime\log}:P^{\log}_{\!U'_\Zar}\to\underline M{}_{|U'_\Zar}$
on some (Zariski) open neighborhood $U'$ of $x$ in $X$ (lemma
\ref{lem_simple-charts-top}(iv.a),(v)). By inspecting the
construction, we find a commutative diagram of monoids :
$$
\xymatrix{ (\tilde u{}^*P^{\log}_{\!U'_\Zar})_\xi \ar[r]
\ar[d]_{(\tilde u{}^*\beta^{\prime\log})_\xi} &
P^{\log}_{\!U_\et,\xi} \ar[d]^{\beta^{\log}_\xi} \\
(\tilde u{}^*\underline M)_\xi \ar@{=}[r] & \underline M{}_{\et,\xi} }
$$
(where $\tilde u{}^*$ is the functor on log structures of
\eqref{subsec_choose-a-top}) whose horizontal arrows and right
vertical arrow are isomorphisms; it follows that
$(\tilde u{}^*\beta^{\prime\log})_\xi$ is an isomorphism as well,
therefore $\tilde u{}^*\beta^{\prime\log}$ restricts to an
isomorphism on some smaller \'etale neighborhood $f:U''\to U'$
of $\xi$ (lemma \ref{lem_simple-charts-top}(iv.c)). Since $f$
is an open morphism, we deduce that the restriction of
$(\tilde u{}^*\beta^{\prime\log})_\xi$ is already an isomorphism on
$f(U'')_\et$, and $f(U'')$ is a (Zariski) open neighborhood of $x$
in $X$. Finally, in light of proposition
\ref{prop_reduce-to-etale}(i), we conclude that the restriction of
$\beta^{\prime\log}$ is an isomorphism on $f(U'')_\Zar$, so
$\underline M$ is coherent, as stated.

Lastly, we suppose that $\underline M{}_\et$ is fine and saturated,
and we wish to show that $\underline M$ is saturated. However, the
assertion can be checked on the stalks, hence let $\xi$ be any
geometric point of $X$; the proof of proposition
\ref{prop_reduce-to-etale} shows that the natural map
$\underline M{}_{|\xi|}^\sharp\to\underline M{}_{\et,\xi}^\sharp$
is bijective, so the assertion follows from lemma
\ref{lem_exc-satura}(ii).
\end{proof}

\begin{corollary} The category $\bcohlog$ admits all finite
limits. More precisely, the limit (in the category of log
structures) of a finite system of coherent log structures,
is coherent.
\end{corollary}
\begin{proof} Let $\Lambda$ be a finite category,
$\cX:=((X_\lambda,\underline{M}_\lambda)~|~\lambda\in\Lambda)$
an inverse system of schemes with coherent log structures
indexed by $\Lambda$. Denote by $Y$ the limit of the system
$(X_\lambda~|~\lambda\in\Lambda)$ of underlying schemes,
and by $\pi_\lambda:Y\to X_\lambda$ the natural morphism,
for every $\lambda\in\Lambda$. It is easily seen that
the limit of $\cX$ is naturally isomorphic to the limit
of the induced system
$\cY:=
((Y,\pi^*_\lambda\underline{M}_\lambda)~|~\lambda\in\Lambda)$,
and in view of lemma \ref{lem_simple-charts-top}(i), we may
therefore replace $\cX$ by $\cY$, and assume that $\cX$
is a finite inverse system in the category $\blog_X$, for
some scheme $X$ (especially, the underlying maps of schemes
are $\one_X$, for every morphism $\lambda\to\mu$ in $\Lambda$).

It suffices then to show that the push-out of two morphisms
$g:\underline{N}\to\underline{M}$,
$h:\underline{N}\to\underline{M}'$ of coherent log structures on
$X$, is coherent. However, notice that the assertion is local on the
site $X_\tau$, hence we may assume that $\underline{N}$ admits a
finite chart $Q_X\to\underline{N}$. Then, thanks to theorem
\ref{th_good-charts}(ii), we may further assume that both $f$ and
$g$ admit finite charts of the form $Q_X\to P_X$ and respectively
$Q_X\to P'_X$, for some finitely generated monoids $P$ and $P'$. We
deduce a natural map $(P\amalg_QP')_X\isom P_X\amalg_{Q_X}P'_X\to
\underline{M}\amalg_{\underline{N}}\underline{M}'$, which is the
sought finite chart.
\end{proof}

\begin{lemma}\label{lem_charts-descend}
In the situation of \eqref{subsec_annoying-general}, suppose that
$X_i$ is quasi-compact for every $i\in I$, and that there exist
$i\in I$, and a coherent log structure $\underline N{}_i$ on $X_i$,
such that the log structure $\underline N:=\pi_i^*\underline N{}_i$
on $X_\tau$ admits a finite chart $\beta:Q_X\to\underline N$. Then
there exist a morphism $\phi:j\to i$ in $I$ and a chart
$\beta_j:Q_{X_j}\to\underline N{}_j:=f_\phi^*\underline N{}_i$ for
$\underline N{}_j$, such that $\pi_j^*\beta_j=\beta$.
\end{lemma}
\begin{proof} After replacing $I$ by $I/i$, we may assume that
$\underline N{}_j$ is well defined for every $j\in I$, and $i$
is the final object of $I$. In this case, notice that
$(X,\underline N)$ is the limit of the cofiltered system of log
schemes $((X_j,\underline N{}_j)~|~j\in I)$. We begin with the
following :

\begin{claim}\label{cl_local-is-good}
In order to prove the lemma, it suffices to show that, for every
$\tau$-point $\xi$ of $X$ there exist $j(\xi)\in I$, a
neighborhood $U_\xi\to X_{j(\xi)}$ of $\pi_{j(\xi)}(\xi)$ in
$X_{j(\xi),\tau}$, and a chart
$\beta_{j(\xi)}:Q_{U_\xi}\to\underline N{}_{j(\xi)|U_\xi}$ such that
$(\one_{U_\xi}\times_{X_{j(\xi)}}\pi_{j(\xi)})^*\beta_{j(\xi)}=\beta$.
\end{claim}
\begin{pfclaim} Clearly we may assume that $U_{\xi}$ is an affine
scheme, for every $\tau$-point $\xi$ of $X$. Under the stated
assumptions, $X$ is quasi-compact, hence we may find a finite set
$\{\xi_1,\dots,\xi_n\}$ of $\tau$-points of $X$, such that
$(U_{\xi_k}\times_{X_{j(\xi_k)}}X\to X~|~k=1,\dots,n)$ is covering
in $X_\tau$. Since $I$ is cofiltered, we may then find $j\in I$ and
morphisms $\phi_k:j\to j(\xi_k)$ for every $k=1,\dots,n$. After
replacing each $\beta_{j(\xi_k)}$ by $f^*_{\phi_k}\beta_{j(\xi_k)}$,
we may assume that $j(\xi_k)=j$ for every $k\leq n$.

In this case, set $\beta_k:=\beta_{j(\xi_k)}$, and $U_k:=U_{\xi_k}$
for every $k\leq n$; let also $U_{kl}:=U_k\times_{X_j}U_l$,
$U^\sim_{kl}:=U_{kl}\times_{X_j}X$ for every $k,l\leq n$, and
denote by $\pi_{kl}:U^\sim_{kl}\to U_{kl}$ the natural projection.
Hence, for every $k,l\leq n$ we have a morphism of $U_{kl}$-monoids :
$$
\beta_{kl}:=\beta_{k|U_{kl}}:Q_{U_{kl}}\to\underline N{}_{j|U_{kl}}
$$
and by construction we have
$\pi^*_{kl}\beta_{kl}=\pi^*_{lk}\beta_{lk}$
under the natural identification : $U_{kl}^\sim\isom U_{lk}^\sim$.
Notice now that $U_{kl}$ is quasi-compact and quasi-separated for
every $k,l\leq n$, hence the natural map
$$
\colim_{i\in I}\Gamma(U_{kl}\times_{X_j}X_i,\underline N{}_i)
\to\Gamma(U_{kl}^\sim,\underline N)
$$
is an isomorphism (lemma \ref{lem_annoying-gen}(ii)).
On the other hand, $\beta_{kl}$ is given by a morphism of
monoids $b_{kl}:Q\to\Gamma(U_{kl},\underline N{}_j)$, and likewise
$\pi^*_{kl}\beta_{kl}$ is given by the morphism
$Q\to\Gamma(U^\sim_{kl},\underline N)$ obtained by composiiton
of $b_{kl}$ and the natural map
$\Gamma(U_{kl},\underline N{}_j)\to\Gamma(U^\sim_{kl},\underline N)$.
By lemma \ref{lem_finite-pres}(ii), we may then find a morphism
$j'\to j$ in $I$ such that the following holds. Set
$$
V_k:=U_k\times_{X_j}X_{j'}
\qquad
V_{kl}:=V_k\times_{X_{j'}}V_l
\qquad
\text{for every $k,l\leq n$}
$$
and let $p_k:V_k\to U_k$ be the projection for every $k\leq n$; let
also $\beta'_k:=p_k^*\beta_k$, which is a chart
$Q_{V_k}\to\underline N{}_{j'|V_k}$ for the restriction of
$\underline N{}_{j'}$. Then
\set\begin{equation}\label{eq_gluing-conditio}
\beta'_{k|V_{kl}}=\beta'_{l|V_{kl}}
\qquad\text{for every $k,l\leq n$}
\end{equation}
under the natural identification $V_{kl}\isom V_{lk}$. By construction,
the system of morphisms $(V_k\times_{X_{j'}}X\to X~|~k=1,\dots,n)$
is covering in $X_\tau$; after replacing $j'$ by a larger index,
we may then assume that the system of morphisms
$(V_k\to X_{j'}~|~k=1,\dots,n)$ is covering in $X_{j',\tau}$
(\cite[Ch.IV, Th.8.10.5]{EGAIV-3}).
In this case, \eqref{eq_gluing-conditio} implies that the
local charts $\beta'_k$ glue to a well defined chart
$\beta_{j'}:Q_{X_{j'}}\to\underline N{}_{j'}$, and a direct
inspection shows that we have indeed $\pi^*_{j'}\beta_{j'}=\beta$.
\end{pfclaim}

Now, denote by $\alpha:Q_X\to\cO_{\!X}$ the composition of
$\beta$ and the structure map of $\underline N$, and let $\xi$
be a $\tau$-point of $X$; we have a natural isomorphism
$$
\underline N{}_\xi=\underline N{}_{i,\pi_i(\xi)}
\otimes_{\underline N{}_{i,\pi_i(\xi)}^\times}\cO_{X,\xi}^\times
\isom\colim_{j\in I}
\underline N{}_{i,\pi_i(\xi)}\otimes_{\underline N{}_{i,\pi_i(\xi)}^\times}
\cO_{\!X_j,\pi_j(\xi_j)}^\times\isom
\colim_{in I}\underline N{}_{j,\pi_j(\xi)}
$$
(\cite[Ch.IV, Prop.18.8.18(ii)]{EGA4}). It follows that $\beta_\xi$
and $\alpha_\xi$ factor through morphisms of monoids
$\beta_{\xi,j}:Q\to\underline N{}_{j,\pi_j(\xi)}$ and
$\alpha_{\xi,j}:Q\to\cO_{\!X_j,\pi_j(\xi)}$ for some $j\in I$ (lemma
\ref{lem_finite-pres}(ii)), and again we may replace $I$ by $I/j$,
after which we may assume that $\beta_{\xi,j}$ and $\alpha_{\xi,j}$
are defined for every $j\in I$. Then $\alpha_{\xi,j}$ extends to a
pre-log structure $\alpha_j:Q_{U_j}\to\cO_{\!U_j}$ on some
neighborhood $U_j\to X_j$ of $\pi_j(\xi)$ in $X_{j,\tau}$ (lemma
\ref{lem_simple-charts-top}(v)), and we may also assume that
$\beta_{\xi,j}$ extends to a morphism of pre-log structures
$\beta_j:(Q_{U_j},\alpha_j)\to\underline N{}_{j|U_j}$ (lemma
\ref{lem_simple-charts-top}(iv.a)). Notice as well that
$$
\underline N{}^\sharp_{j,\xi}\simeq\underline N{}_\xi^\sharp
\quad\text{and}\quad
\beta_{\xi,j}^{-1}\underline N{}^\times_{j,\pi_j(\xi)}=
\beta_\xi^{-1}\underline N{}_\xi^\times
\qquad
\text{for every $j\in I$}
$$
and since $\beta_\xi$ induces an isomorphism
$Q/\beta_\xi^{-1}\underline N{}^\times_\xi\isom
\underline N{}^\sharp_\xi$, we deduce that $\beta_{j,\xi}$
(which is the same as $\beta_{\xi,j}$) induces an isomorphism
$Q/\alpha_{j,\xi}^{-1}\cO^\times_{\!X_j,\pi_j(\xi)}\isom
\underline N{}^\sharp_{j,\xi}$ for every $j\in I$. In turn,
it then follows from lemma \ref{lem_check-iso} that $\beta_{j,\xi}$
induces an isomorphism
$(Q_{U_j},\alpha_j)^{\log}_{\pi_j(\xi)}\isom
\underline N{}_{j,\pi_j(\xi)}$.

Next, by lemma \ref{lem_simple-charts-top}(iv.b,c) we may find, for
every $j\in I$, a neighborhood $U'_j\to U_j$ of $\xi$ in $X_\tau$,
such that the restriction
$(Q_{U'_j},\alpha_{j|U'_j})\to\underline N{}_{j|U'_j}$ of
$\beta_j$ is a chart for $\underline N{}_{j|U'_j}$.

By construction, the morphism
$((\one_{U_j}\times_{X_j}\pi_j)^*\beta_j)_\xi:Q\to\underline N{}_\xi$
is the same as $\beta_\xi$, so by lemma
\ref{lem_simple-charts-top}(iv.b) we may find a neighborhood
$V_j\to U'_j\times_{X_j}X$ of $\xi$ in $X_\tau$, such that
$$
((\one_{U_j}\times_{X_j}\pi_j)^*\beta_j)_{|V_j}=\beta_{|V_j}
\qquad\text{and}\qquad
((\one_{U_j}\times_{X_j}\pi_j)^*\alpha_j)_{|V_j}=\alpha_{|V_j}.
$$
\begin{claim}\label{cl_combination}
In the situation of \eqref{subsec_annoying-general}, let $Y\to X$ be
an object of the site $X_\tau$, with $Y$ quasi-compact and
quasi-separated. We have :
\begin{enumerate}
\item
There exist $i\in I$, an object $Y_i\to X_i$ of $X_{i,\tau}$, and
an isomorphism of $X$-schemes $Y\isom Y_i\times_{X_i}X$.
\item
Moreover, if $Y\to X$ is covering in $X_\tau$, then we may find
$i\in I$ and $Y_i\to X_i$ as in (i), which is covering in
$X_{i,\tau}$.
\end{enumerate}
\end{claim}
\begin{pfclaim}(i) is obtained by combining
\cite[Ch.IV, Th.8.8.2(ii)]{EGAIV-3} and \cite[Ch.IV,
Prop.17.7.8(ii)]{EGA4} (and \cite[Ch.IV, Cor.8.6.4]{EGAIV-3} if
$\tau=\Zar$). Assertion (ii) follows from (i) and \cite[Ch.IV,
Th.8.10.5]{EGAIV-3}.
\end{pfclaim}

By claim \ref{cl_combination}(i), we may assume that
$V_j=U''_j\times_{X_j}X$ for some neighborhood $U''_j\to U'_j$ of
$\pi_j(\xi)$ in $X_{j,\tau}$. To conclude, it suffices to invoke
claim \ref{cl_local-is-good}.
\end{proof}

\begin{proposition}\label{prop_2-colim-for-logs}
In the situation of \eqref{subsec_annoying-general}, suppose that
$X_i$ is quasi-compact for every $i\in I$. Then the natural functor :
\set\begin{equation}\label{eq_this_functor}
\Pscolim{I}\bcohlog_{X_i}\to\bcohlog_{X}
\end{equation}
is an equivalence.
\end{proposition}
\begin{proof} To begin with, we show :
\begin{claim}\label{cl_it-is-faithful}
The functor \eqref{eq_this_functor} is faithful. Namely, for a given
$i\in I$, let $\underline M{}_i$ and $\underline N{}_i$ be two coherent
log structures on $X_i$, and $f_i,g_i:\underline M{}_i\to\underline N{}_i$
two morphisms, such that $\pi_i^*f_i$ agrees with $\pi_i^*g_i$.
Then there exists a morphism $\psi:j\to i$ in $I$ such that
$f^*_\psi f_i=f^*_\psi g_i$.
\end{claim}
\begin{pfclaim} Set $\underline M:=\pi_i^*\underline M{}_i$, and define
likewise the coherent log structure $\underline N$ on $X$. For any
$\tau$-point $\xi$ of $X$, pick a neighborhood $U_\xi\to X_i$ of
$\pi_i(\xi)$ in $X_i$, and a finite chart
$\beta:P_{U_\xi}\to\underline M{}_{i|U_\xi}$ for the restriction of
$\underline M{}_i$. The morphisms $f_{i|U_\xi}$ and $g_{i|U_\xi}$ are
determined by the induced maps $\phi:=\Gamma(U_\xi,f_i)\circ\beta$
and $\gamma:=\Gamma(U_\xi,g_i)\circ\beta$, and the assumption means
that the composition of $\phi$ and the natural map
$\Gamma(U_\xi,\underline N{}_i)\to\Gamma(U_\xi\times_{X_i}X,\underline N)$
equals the composition of $\gamma$ with the same map.

It then follows from lemmata \ref{lem_annoying-gen}(ii) and
\ref{lem_finite-pres}(ii) that there exists a morphism
$\psi:i(\xi)\to i$ in $I$ such that the composition of $\phi$ with
the natural map $\Gamma(U_\xi,\underline N{}_i)\to
\Gamma(U_\xi\times_{X_i}X_{i(\xi)},f^*_\psi\underline N)$
equals the composition of $\gamma$ with the same map; in other
words, if we set $U'_\xi:=U_\xi\times_{X_i}X_{i(\xi)}$, we have
$f^*_\psi f_{i|U'_\xi}=f^*_\psi g_{i|U'_\xi}$. Next, since $X$ is
quasi-compact, we may find finitely many $\tau$-points
$\xi_1,\dots,\xi_n$ such that the family
$(U'_{\xi_k}\times_{X_{i(\xi_k)}}X\to X)$ is covering in $X_\tau$.
Then, by \cite[Ch.IV, Th.8.10.5]{EGAIV-3} we may find $j\in I$ and
morphisms $\psi_k:j\to i(\xi_k)$ in $I$, for $k=1,\dots,n$, such
that the induced  family
$(U''_k:=U'_{\xi_k}\times_{X_{i(\xi_k)}}X_j\to X_j)$ is covering in
$X_{j,\tau}$. By construction we have
$$
f^*_{\psi_k\circ\psi}f_{i|U''_k}=f^*_{\psi_k\circ\psi}g_{i|U''_k}
\qquad\text{for $k=1,\dots,n$}
$$
therefore $f^*_{\psi_k\circ\psi}f_i=f^*_{\psi_k\circ\psi}g_i$, as
required.
\end{pfclaim}

\begin{claim}\label{cl_coralie}
The functor \eqref{eq_this_functor} is full. Namely, let $i\in I$
and $\underline M_i$, $\underline N_i$ as in claim
\ref{cl_it-is-faithful}, and suppose that $f:\pi_i^*\underline M{}_i
\to\pi_i^*\underline N{}_i$ is a given morphism of log structures;
then there exist a morphism $\psi:j\to i$ in $I$, and a morphism of
log structures $f_j:f^*_\psi\underline M_i\to f^*_\psi\underline N{}_i$
such that $\pi_j^*f_j=f$.
\end{claim}
\begin{pfclaim} Indeed, by theorem \ref{th_good-charts}(i) we may find
a covering family $(U_\lambda\to X~|~\lambda\in\Lambda)$ in
$X_\tau$, and for each $\lambda\in\Lambda$ a finite chart
$$
\beta_\lambda:P_{\lambda,U_\lambda}\to\underline M{}_{|U_\lambda}
\quad
\gamma_\lambda:Q_{\lambda,U_\lambda}\to\underline N{}_{|U_\lambda}
\quad
\theta_\lambda:P_\lambda\to Q_\lambda
$$
for the restriction $f_{|U_\lambda}$. Clearly we may assume that
each $U_\lambda$ is affine, and since $X$ is quasi-compact, we may
assume as well that $\Lambda$ is a finite set; in this case, claim
\ref{cl_combination} implies that there exist a morphism $j\to i$
in $I$, a covering family $(U_{j,\lambda}\to
X_j~|~\lambda\in\Lambda)$, and isomorphism of $X$-schemes
$U_\lambda\isom U_{j,\lambda}\times_{X_j}X$ for every
$\lambda\in\Lambda$. Next, lemma \ref{lem_charts-descend} says that,
after replacing $j$ by some larger index, we may assume that for
every $\lambda\in\Lambda$ there exist charts
$$
\beta_{j,\lambda}:P_{\lambda,U_{j,\lambda}}\to\underline M{}_j:=
f^*_\psi\underline M{}_i
\qquad\text{and}\qquad
\gamma_{j,\lambda}:Q_{\lambda,U_{j,\lambda}}\to\underline N{}_j:=
f^*_\psi\underline N{}_i
$$
with $\pi_j^*\beta_{j,\lambda}=\beta_\lambda$ and
$\pi_j^*\gamma_{j,\lambda}=\gamma_\lambda$. Set
$f_{j,\lambda}:=\theta_{\lambda,U_{j,\lambda}}^{\log}:
\underline M{}_j\to\underline N{}_j$; by construction we have :
$$
(\one_{U_{j,\lambda}}\times_{X_j}\pi_j)^*f_{j,\lambda}=f_{|U_\lambda}
\qquad \text{for every $\lambda\in\Lambda$}.
$$
Next, for every $\lambda,\mu\in\Lambda$, let
$U_{j,\lambda\mu}:=U_{j,\lambda}\times_{X_j}U_{j,\mu}$; we deduce,
for every $\lambda,\mu\in\Lambda$, two morphisms of log structures
$f_{j,\lambda|U_{j,\lambda\mu}},f_{j,\mu|U_{\lambda\mu}}:
\underline M{}_{j|U_{\lambda\mu}}\to\underline N{}_{j|U_{\lambda\mu}}$,
which agree after pull back to $U_{j,\lambda\mu}\times_{X_j}X$.
By applying claim \ref{cl_it-is-faithful} to the cofiltered
system of schemes
$(U_{j,\lambda\mu}\times_{X_j}X_{j'}~|~j'\in I/j)$, we may then
achieve -- after replacing $j$ by some larger index -- that
$f_{j,\lambda|U_{j,\lambda\mu}}=f_{j,\mu|U_{\lambda\mu}}$, in which
case the system $(f_{j,\lambda}~|~\lambda\in\Lambda)$ glues to a
well defined morphism $f_j$ as sought.
\end{pfclaim}

Finally, let us show that \eqref{eq_this_functor} is essentially
surjective. Indeed, let $\underline M$ be a coherent log structure
on $X$; let us pick a covering family
$\cU:=(U_\lambda\to X~|~\lambda\in\Lambda)$ and finite charts
$\beta_\lambda:P_{\lambda,U_\lambda}\to\underline M{}_{|U_\lambda}$
for every $\lambda\in\Lambda$. As in the foregoing, we may assume
that each $U_\lambda$ is affine, and $\Lambda$ is a finite set, in
which case, according to claim \ref{cl_combination}(ii) we may find
$i\in I$ and a covering family
$(U_{i,\lambda}\to X_i~|~\lambda\in\Lambda)$ with isomorphisms
of $X$-schemes $U_{i,\lambda}\times_{X_i}X\isom U_\lambda$ for
every $\lambda\in\Lambda$; for every $j\in I/i$ and every
$\lambda\in\Lambda$, let us set
$U_{j,\lambda}:=U_{i,\lambda}\times_{X_i}X_j$. The composition
$\alpha_\lambda:P_{\lambda,U_\lambda}\to\cO_{U_\lambda}$ of
$\beta_\lambda$ and the structure map of $\underline M{}_{|U_\lambda}$
is determined by a morphism of monoids
$$
P_\lambda\to\Gamma(U_\lambda,\cO_{U_\lambda})=
\colim_{i\in I/i}\Gamma(U_{j,\lambda},\cO_{U_{j,\lambda}}).
$$
Then, as usual, lemma \ref{lem_finite-pres}(ii) implies that there
exists $j\in I/i$ such that $\alpha_\lambda$ descends to a pre-log
structure $P_{\lambda,U_{j,\lambda}}\to\cO_{U_{j,\lambda}}$ on
$U_{j,\lambda}$, whose associated log structure we denote by
$\underline M{}_{j,\lambda}$. For every $\lambda,\mu\in\Lambda$,
let $U_{j,\lambda\mu}:=U_{j,\lambda}\times_{X_j}U_{j,\mu}$; by
construction we have isomorphisms
\set\begin{equation}\label{eq_natasha}
(\one_{U_{j,\lambda\mu}}\times_{X_j}\pi_j)^*
\underline M{}_{j,\lambda|U_{j,\lambda\mu}}\isom
(\one_{U_{j,\mu\lambda}}\times_{X_j}\pi_j)^*
\underline M{}_{j,\mu|U_{j,\mu\lambda}}
\qquad
\text{for every $\lambda,\mu\in\Lambda$}.
\end{equation}
By applying claim \ref{cl_coralie} to the cofiltered system
$(U_{j,\lambda\mu}\times_{X_j}X_{j'}~|~j'\in I/j)$, we can then
obtain -- after replacing $j$ by a larger index -- isomorphisms
$\omega_{\lambda\mu}:\underline M{}_{j,\lambda|U_{j,\lambda\mu}}\isom
M_{j,\mu|U_{j,\mu\lambda}}$ for every $\lambda,\mu\in\Lambda$, whose
pull-back to $U_{j,\lambda\mu}\times_{X_j}X$ are the isomorphisms
\eqref{eq_natasha}. Lastly, in view of claim
\ref{cl_it-is-faithful}, we may achieve -- after further replacement
of $j$ by a larger index -- that the system $\cD:=
(\underline M{}_{j,\lambda},\omega_{\lambda\mu}~|~\lambda,\mu\in\Lambda)$
is a descent datum for the fibration $F$ of \eqref{subsec_special-schs},
whose pull-back to $X$ is isomorphic to the natural descent datum
for $\underline M$, associated to the covering family $\cU$. Then
$\cD$ glues to a log structure $\underline M{}_j$, such that
$\pi_j^*\underline M{}_j\simeq\underline M$.
\end{proof}

\begin{corollary}\label{cor_notime-now}
In the situation of \eqref{subsec_lambdasmus}, suppose that
$X_\lambda$ is quasi-compact for every $\lambda\in\Ob(\Lambda)$,
and let
$$
(g,\log g):(Y,\underline N)\to(X,\underline M)
$$
be a morphism of log schemes with coherent log structures.
Then there exist $\lambda\in\Lambda$, and a morphism
$$
(g_\lambda,\log g_\lambda):(Y_\lambda,\underline N{}_\lambda)
\to(X_\lambda,\underline M{}_\lambda)
$$
of log schemes with coherent log structures, such that
$\log g=\psi_\lambda^*\log g_\lambda$.
\end{corollary}
\begin{proof} The assertion is an immediate
consequence of proposition \ref{prop_2-colim-for-logs}.
\end{proof}

\begin{corollary}\label{cor_descend-chart-from-infty}
In the situation of \eqref{subsec_lambdasmus}, let
$0\in\Ob(\Lambda)$ be an index such that
$\underline N_0$ and $\underline M_0$ are coherent log
structures on $Y_0$, and respectively $X_0$, and
$(g_0,\log g_0):(Y_0,\underline N{}_0)\to(X_0,\underline M{}_0)$
a morphism of log schemes. Suppose also that $X_\lambda$ is
quasi-compact (as well as quasi-separated) for every
$\lambda\in\Ob(\Lambda)$; then we have :
\begin{enumerate}
\item
If the morphism of log schemes
$$
(g,\psi_0^*\log g_0):Y\times_{Y_0}(Y_0,\underline N{}_0)
\to X\times_{X_0}(X_0,\underline M{}_0)
$$
admits a chart $(\omega,\omega',\theta:P\to Q)$,
there exists a morphism $u:\lambda\to 0$ in $\Lambda$ such that
the morphism of log schemes
$$
(g_\lambda,\psi_u^*\log g_0):Y_\lambda\times_{Y_0}(Y_0,\underline N{}_0)
\to X_\lambda\times_{X_0}(X_0,\underline M{}_0)
$$
admits a chart $(\omega_\lambda,\omega'_\lambda,\theta)$.
\item
If $(g,\psi_0^*\log g_0)$ is a log flat (resp. saturated)
morphism of fine log schemes, there exists a morphism
$u:\lambda\to 0$ in $\Lambda$ such that $(g_\lambda,\psi_u^*\log g_0)$
is a log flat (resp. saturated) morphism of fine log schemes.
\end{enumerate}
\end{corollary}
\begin{proof}(i): Under the current assumptions, $X$ and
$Y$ are quasi-compact. In view of lemma \ref{lem_charts-descend},
we may then find a morphism $v:\mu\to 0$, such that $\omega$
and $\omega'$ descend to charts
$\omega_v:P_{X_\mu}\to\phi_v^*\underline M{}_0$ and
$\omega'_v:Q_{Y_\mu}\to\psi_v^*\underline N{}_0$. We deduce two
morphisms of pre-log structures
$P_{Y_\mu}\to\psi^*_v\underline N{}_0$, namely
$$
\beta_1:=\psi^*_v(\log g_0)\circ g_\mu^*\omega_v \qquad \text{and}
\qquad \beta_2:=\omega'_v\circ\theta_{Y_\mu}
$$
and by construction we have
$\psi_\mu^*\beta_1^{\log}=\psi_\mu^*\beta_2^{\log}$. By proposition
\ref{prop_2-colim-for-logs} it follows that there exists
$w:\lambda\to\mu$ such that
$\psi_w^*\beta_1^{\log}=\psi_w^*\beta_2^{\log}$. The latter means
that $(\phi^*_w\omega_v,\psi^*_w\omega'_v,\theta)$ is a chart for
the morphism of log schemes
$(g_\lambda,\psi_{v\circ w}^*\log g_0):
(Y_\lambda,\psi_{v\circ w}^*\underline N{}_0)\to
(X_\lambda,\phi_{v\circ w}^*\underline M{}_0)$, {\em i.e.} the
claim holds with $u:=v\circ w$.

(ii): Suppose therefore that $(g,\psi_0^*\log g_0)$ is a
log flat (resp.saturated) morphism, and let
$\cU:=(U_i\to X~|~i\in I)$ be a covering family for
$X_\tau$, such that $U_i\times_{X_0}(X_0,\underline M{}_0)$
admits a finite (resp. fine) chart, and $U_i$ is affine for
every $i\in I$.
Since $X$ is quasi-compact, we may assume that $I$ is a
finite set, and then there exists $\lambda\in\Lambda$ such that
$\cU$ descends to a covering family
$\cU_\lambda:=(U_{\lambda,i}\to X_\lambda~|~i\in I)$ for
$X_{\lambda,\tau}$ (claim \ref{cl_combination}(ii)). After
replacing $\Lambda$ by $\Lambda/\lambda$, we may assume that
$\lambda=0$, in which case we set
$U_{\lambda,i}:=U_{0,i}\times_{X_0}X_\lambda$ for every object
$\lambda$ of $\Lambda$, and every $i\in I$. Clearly, it suffices
to show that there exists $u:\lambda\to 0$ such that
$U_{\lambda,i}\times_{X_\lambda}(g_\lambda,\psi_u^*\log g_0)$
is flat (resp. saturated) for every $i\in I$. Set
$Y'_{\lambda,i}:=Y_\lambda\times_{X_\lambda}U_{\lambda,i}$ for every
$\lambda\in\Lambda$; we may then replace the cofiltered
system $\underline X$ and $\underline Y$, by respectively
$(U_{\lambda,i}~|~\lambda\in\Lambda)$ and
$(Y'_{\lambda,i}~|~\lambda\in\Lambda)$, which allows to assume
from start, that $X\times_{X_0}(X_0,\underline M{}_0)$
admits a finite (resp. fine) chart. In this case, lemma
\ref{lem_charts-descend} allows to further reduce to the case
where $\underline M{}_0$ admits a finite (resp. fine) chart.

Then, by theorem \ref{th_good-charts}(iii), we may
find a covering family
$\cV:=(V_j\to Y~|~j\in J)$ for $Y_\tau$, consisting
of finitely many affine schemes $V_j$, and for every $j\in J$,
a flat (resp. saturated) and fine chart for the induced
morphism $V_j\times_{Y_0}(Y_0,\underline N{}_0)\to
X\times_{X_0}(X_0,\underline M{}_0)$. As in the foregoing,
after replacing $\Lambda$ by some category $\Lambda/\lambda$,
we may assume that $\cV$ descends to a covering family
$\cV_0:=(V_{0,j}\to Y_0~|~j\in J)$ for $Y_{0,\tau}$, in which
case we set $V_{\lambda,j}:=V_{0,j}\times_{Y_0}Y_\lambda$
for every $\lambda\in\Lambda$.
Clearly, it suffices to show that there exists
$\lambda\in\Lambda$ such that the induced morphism
$V_{\lambda,j}\times_{Y_0}(Y_0,\underline N{}_0)\to
X_\lambda\times_{X_0}(X_0,\underline M{}_0)$ is log flat
(resp. saturated). Thus, we may replace $\underline Y$ by
the cofiltered system $(V_{\lambda,j}~|~\lambda\in\Lambda)$,
which allows to assume that $(g,\psi_0^*\log g_0)$
admits a flat (resp. saturated) and fine chart. In this case,
the assertion follows from (i).
\end{proof}

\begin{proposition}\label{prop_quasi-fine}
The inclusion functors :
$$
\bqflog\to\bqcohlog \qquad \bqfslog\to\bqflog
$$
admit right adjoints :
$$
\bqcohlog\to\bqflog\ :\
(X,\underline M)\mapsto(X,\underline M)^\mathrm{qf}
\qquad
\bqflog\to\bqfslog\ :\
(X,\underline M)\mapsto(X,\underline M)^\mathrm{qfs}.
$$
\end{proposition}
\begin{proof} Let $(X,\underline M)$ be a scheme with quasi-coherent
(resp. quasi-fine) log structure. We need to construct a morphism of log
schemes
$$
\phi:(X,\underline M)^\mathrm{qf}\to(X,\underline M)
\qquad
\text{(resp. $\phi:(X,\underline M)^\mathrm{qfs}\to(X,\underline M)$)}
$$
such that $(X,\underline M)^\mathrm{qf}$ (resp. $(X,\underline
M)^\mathrm{qfs}$) is a quasi-fine (resp. qfs) log scheme, and the
following holds. Every morphism of log schemes $\psi:(Y,\underline
N)\to(X,\underline M)$ with $(Y,\underline N)$ quasi-fine (resp.
qfs), factors uniquely through $\phi$. To this aim, suppose first
that $\underline M$ admits a chart (resp. a quasi-fine chart)
$\alpha:P_X\to\underline M$. By lemma \ref{lem_time}(i),
$\alpha$ determines an isomorphism
$$
(X,\underline M)\isom\Spec(\Z,P)\times_{\Spec\,\Z[P]}X.
$$
Since $\underline N$ is integral (resp. and saturated), the
morphism $(Y,\underline N)\to\Spec(\Z,P)$ induced by $\psi$
factors uniquely through $\Spec(\Z,P^\intg)$ (resp.
$\Spec(\Z,P^\sat)$) (lemma \ref{lem_time}(ii)). Taking into
account lemma \ref{lem_simple-charts-top}(iii), it follows
easily that we may take
$$
(X,\underline M)^\mathrm{qf}:=
\Spec(\Z,P^\intg)\times_{\Z[P]}X
\qquad
(X,\underline M)^\mathrm{qfs}:=
\Spec(\Z,P^\sat)\times_{\Z[P]}X
$$
Next, notice that the universal property of
$(X',\underline M'):=(X,\underline M)^\mathrm{qf}$
(resp. of $(X',\underline M'):=(X,\underline M)^\mathrm{qfs}$)
is local on $X_\tau$ : namely, suppose that $(X',\underline M')$
has already been found, and let $U\to X$ be an object of $X_\tau$,
with a morphism $(Y,\underline N)\to(U,\underline M_{|U})$ from a
quasi-fine (resp. qfs) log scheme; there follows a unique morphism
$$
(Y,\underline N)\to
(U,\underline M_{|U})\times_{(X,\underline M)}(X',\underline M')
\isom U\times_X(X',\underline M')
$$
(notation of \eqref{eq_prod-sch-log-sch}). Thus
$(U,\underline M_{|U})^\mathrm{qf}\simeq
U\times_X(X,\underline M)^\mathrm{qf}$ (resp.
$(U,\underline M_{|U})^\mathrm{qfs}\simeq
U\times_X(X,\underline M)^\mathrm{qfs}$,)
since the latter is a quasi-fine (resp. qfs) log scheme. Therefore,
for a general quasi-coherent (resp. quasi-fine) log structure
$\underline M$, choose a covering family
$(U_\lambda\to X~|~\lambda\in\Lambda)$ such that
$(U_\lambda,\underline M_{|U_\lambda})$ admits a chart
for every $\lambda\in\Lambda$; it follows that the family
$\cU:=
((U_\lambda,\underline M_{|U_\lambda})^\mathrm{qf}~|~\lambda\in\Lambda)$,
together with the natural isomorphisms :
$$
U_\mu\times_X(U_\lambda,\underline M_{|U_\lambda})^\mathrm{qf}\simeq
U_\lambda\times_X(U_\mu,\underline M_{|U_\mu})^\mathrm{qf}
\qquad\text{for every $\lambda,\mu\in\Lambda$}
$$
is a descent datum for the fibred category over $X_\tau$,
whose fibre over any object $U\to X$ is the category of
affine $U$-schemes endowed with a log structure
(resp. likewise for the family
$\cU:=
((U_\lambda,\underline M_{|U_\lambda})^\mathrm{qfs}~|~\lambda\in\Lambda)$).
Using faithfully flat descent (\cite[Exp.VIII, Th.2.1]{SGA1}),
one sees that $\cU$ comes from a quasi-fine (resp. qfs) log scheme
which enjoys the sought universal property.
\end{proof}

\begin{remark}\label{rem_quasi-fine}
(i)\ \ 
By inspecting the proof of proposition \ref{prop_quasi-fine},
we see that the quasi-fine log scheme associated to a coherent
log scheme, is actually fine, and the qfs log scheme associated
to a fine log scheme, is a fs log scheme (one applies corollary
\ref{cor_fragment-Gordon}(i)). Hence we obtain functors
$$
\bcohlog\to\bflog
\ :\ (X,\underline M)\mapsto(X,\underline M)^\mathrm{f}
\qquad
\bflog\to\bfslog\ :\
(X,\underline M)\mapsto(X,\underline M)^\mathrm{fs}
$$
which are right adjoint to the inclusion functors
$\bflog\to\bcohlog$ and $\bfslog\to\bflog$.

(ii)\ \ 
Notice as well that, for every log scheme $(X,\underline M)$
with quasi-coherent (resp. fine) log structure, the morphism
of schemes underlying the counit of adjunction
$(X,\underline M)^\mathrm{qf}\to(X,\underline M)$ (resp.
$(X,\underline M)^\mathrm{fs}\to(X,\underline M)$) is a closed
immersion (resp. is finite).

(iii)\ \ 
Furthermore, let $f:Y\to X$ be any morphism of schemes; the
proof of proposition \ref{prop_quasi-fine} also yields a
natural isomorphism of $(Y,f^*\underline M)$-schemes :
$$
(Y,f^*\underline M)^\mathrm{qf}\isom
Y\times_X(X,\underline M)^\mathrm{qf}
\qquad
\text{(resp. $(Y,f^*\underline M)^\mathrm{qfs}\isom
Y\times_X(X,\underline M)^\mathrm{qfs}$)}.
$$

(iv)\ \
Let $(X,\underline M)$ be a quasi-fine log scheme, and suppose that
$X$ is a normal, irreducible scheme, and $(X,\underline M)_\tr$
is a dense subset of $X$. Denote by $X^\mathrm{qfs}$ the scheme
underlying $(X,\underline M)^\mathrm{qfs}$; then we claim
that the projection $X^\mathrm{qfs}\to X$ (underlying the counit
of adjunction) admits a natural section :
$$
\sigma_X:X\to X^\mathrm{qfs}.
$$
The naturality means that if $f:(X,\underline M)\to(Y,\underline N)$
is any morphism of quasi-fine log schemes, where $Y$ is also normal
and irreducible, and $(Y,\underline N)_\tr$ is dense in $Y$,
then the induced diagram of schemes :
\set\begin{equation}\label{eq_section-natu}
{\diagram
X \ar[r]^-{\sigma_X} \ar[d]_f & X^\mathrm{qfs} \ar[d]^{f^\mathrm{qfs}} \\
Y \ar[r]^-{\sigma_Y} & Y^\mathrm{qfs}
\enddiagram}\end{equation}
commutes, and therefore it is cartesian, by virtue of (iii).
Indeed, suppose first that $X$ is affine, say $X=\Spec\,A$ for
some normal domain $A$, and $\underline M$ admits an integral
chart, given by a morphism $\beta:P\to A$, for some integral
monoid $P$; we have to exhibit a ring homomorphism
$P^\sat\otimes_PA\to A$, whose composition with the natural
map $A\to P^\sat\otimes_PA$ is the identity of $A$. The latter
is the same as the datum of a morphism of monoids $P^\sat\to A$
whose restriction to $P$ agrees with $\beta$. However, since
the trivial locus of $\underline M$ is dense in $X$, the image
of $P$ in $A$ does not contain $0$, hence $\beta$ extends to a
group homomorphism $\beta^\gp:P^\gp\to\Frac(A)^\times$; since
$A$ is integrally closed in $\Frac(A)$, we have
$\beta^\gp(P^\sat)\subset A$, as required. Next, suppose that
$U\to X$ is an object of $X_\tau$, with $U$ also affine and
irreducible; then $U$ is normal and the trivial locus of
$(U,\underline M{}_{|U})$ is dense in $U$. Thus, the foregoing
applies to $(U,\underline M_{|U})$ as well, and by inspecting
the constructions we deduce a natural identification :
$$
\sigma_U=\one_U\times_X\sigma_X.
$$
Lastly, for a general $(X,\underline M)$, we can find a covering
family $(U_\lambda\to X~|~\lambda\to\Lambda)$ in $X_\tau$, such
that $U_\lambda$ is affine, and
$(U_\lambda,\underline M{}_{|U_\lambda})$ admits an integral
chart for every $\lambda\in\Lambda$; proceeding as above, we
obtain a system of morphisms
$(\sigma_\lambda:U_\lambda\to U^\mathrm{qfs}~|~\lambda\in\Lambda)$,
as well as natural identifications :
$$
\one_{U_\mu}\times_X\sigma_\lambda=
\one_{U_\lambda}\times_X\sigma_\mu
\qquad
\text{for every $\lambda,\mu\in\Lambda$}.
$$
In other words, we have defined a descent datum for the category
fibred over $X_\tau$, whose fibre over any object $U\to X$ is
the category of all morphisms of schemes $U\to U^\mathrm{qfs}$.
By invoking faithfully flat descent (\cite[Exp.VIII, Th.2.1]{SGA1}),
we see that this descent datum yields a morphism
$\sigma_X:X\to X^\mathrm{qfs}$ such that
$\one_{U_\lambda}\times_X\sigma_X=\sigma_\lambda$ for every
$\lambda\in\Lambda$. The verification that $\sigma_X$ is a
section of the projection $X^\mathrm{qfs}\to X$, and that
\eqref{eq_section-natu} commutes, can be carried out locally
on $X_\tau$, in which case we can assume that $\underline M$
admits a chart as above, and one can check explicitly these
assertions, by inspecting the constructions.
\end{remark}

\subsection{Logarithmic differentials and smooth morphisms}
\label{sec_smooth-log-sch}
In this section we introduce the logarithmic version of
the usual sheaves of relative differentials, and we study
some special classes of morphisms of log schemes.

\begin{definition}
Let $(X,\underline{M}\xrightarrow{\alpha}\cO_{\!X})$ and
$(Y,\underline{N}\xrightarrow{\beta}\cO_{\!Y})$ be two
schemes with pre-log structures, and
$f:(X,\underline{M})\to(Y,\underline{N})$ a morphism of
schemes with pre-log structures. Let also $\cF$ be an
$\cO_{\!X}$-module. An {\em $f$-linear derivation of
$\underline M$ with values in $\cF$} is a pair
$(\partial,\log\partial)$ consisting of maps of sheaves :
$$
\partial:\cO_{\!X}\to\cF
\qquad
\log\partial:\log\underline M\to\cF
$$
such that :
\begin{itemize}
\item
$\partial$ is a derivation (in the usual sense).
\item
$\log\partial$ is a morphism of sheaves of (additive) monoids
on $X_\tau$.
\item
$\partial\circ f^\natural=0$ and $\log\partial\circ\log f=0$.
\item
$\partial\circ\alpha(m)=\alpha(m)\cdot\log\partial(m)$
for every object $U$ of $X_\tau$, and every $m\in\underline M(U)$.
\end{itemize}
The set of all $f$-linear derivations with values
in $f$ shall be denoted by :
$$
\Der_{(Y,\underline N)}((X,\underline M),\cF).
$$
The $f$-linear derivations shall also be called {\em
$(Y,\underline N)$-linear derivations}, when there is no danger
of ambiguity.
\end{definition}

\sset\subsubsection{}
The set $\Der_{(Y,\underline N)}((X,\underline M),\cF)$ is
clearly functorial in $\cF$, and moreover, for any object
$U$ of $X_\tau$, any $s\in\cO_{\!X}(U)$, and any $f$-linear
derivation $(\partial,\log\partial)$, the restriction
$(s\cdot\partial_{|U},s\cdot\log\partial_{|U})$ is an
element of $\Der_{(Y,\underline N)}((U,\underline M_{|U}),\cF)$,
hence the rule
$U\mapsto\Der_{(Y,\underline N)}((U,\underline M_{|U}),\cF)$
defines an $\cO_{\!X}$-module :
$$
\cDer_{(Y,\underline N)}((X,\underline M),\cF).
$$
In case $\underline M=\cO^\times_{\!X}$ and
$\underline N=\cO^\times_{\!Y}$ are the trivial log structures,
the $f$-linear derivations of $\underline M$ are the same as
the usual $f$-linear derivations, {\em i.e.} the natural map
\set\begin{equation}\label{eq_triv-deriv}
\cDer_{(Y,\cO^\times_{\!Y})}((X,\cO^\times_{\!X}),\cF)\to
\cDer_{Y}(X,\cF)
\end{equation}
is an isomorphism.
In the category of usual schemes, the functor of derivations
is represented by the module of relative differentials.
This construction extends to the category of schemes with
pre-log structures. Namely, let us make the following :

\begin{definition}\label{def_log-Kahler}
Let $(X,\underline M\xrightarrow{\alpha}\cO_{\!X})$ and
$(Y,\underline N\xrightarrow{\beta}\cO_{\!Y})$ be two schemes with
pre-log structures, and $f:(X,\underline M)\to(Y,\underline N)$
a morphism of schemes with pre-log structures. The
{\em sheaf of logarithmic differentials\/} of $f$ is the
$\cO_{\!X}$-module :
$$
\Omega^1_{X/Y}(\log\underline M/\underline N):=
(\Omega^1_{X/Y}\oplus(\cO_X\otimes_\Z\log\underline M^\gp))/R
$$
where $R$ is the $\cO_{\!X}$-submodule generated locally
on $X_\tau$ by local sections of the form :
\begin{itemize}
\item
$(d\alpha(a),-\alpha(a)\otimes\log a)$ with $a\in\underline M(U)$
\item
$(0,1\otimes\log a)$ with
$a\in\Img((f^{-1}\underline N)(U)\to\underline M(U))$
\end{itemize}
where $U$ ranges over all the objects of $X_\tau$
(here we use the notation of \eqref{sec_toric}).
For arbitrary $a\in\underline M(U)$, the class of
$(0,1\otimes\log a)$ in
$\Omega^1_{X/Y}(\log\underline M/\underline N)$
shall be denoted by $d\log a$.
\end{definition}

\sset\subsubsection{}
It is easy to verify that
$\Omega^1_{X/Y}(\log\underline M/\underline N)$
represents the functor
$$
\cF\mapsto\Der_{(Y,\underline N)}((X,\underline M),\cF)
$$
on $\cO_{\!X}$-modules. Consequently, \eqref{eq_triv-deriv}
translates as a natural isomorphism of $\cO_{\!X}$-modules :
\set\begin{equation}\label{eq_triv-differ}
\Omega^1_{X/Y}\isom\Omega^1_{X/Y}(\log\cO^\times_{\!X}/\cO^\times_{\!Y}).
\end{equation}
Furthermore, let us fix a scheme with pre-log structure
$(S,\underline N)$, and define the category :
$$
\bprelog/(S,\underline N)
$$
as in \eqref{subsec_opposing}.
Also, let $\logMod/(S,\underline N)$ be the category whose
objects are all the pairs $((X,\underline M),\cF)$, where
$(X,\underline M)$ is a $(S,\underline N)$-scheme, and
$\cF$ is an $\cO_{\!X}$-module. The morphisms
$$
((X,\underline M),\cF)\to((Y,\underline N),\cG)
$$
in $\logMod/(S,\underline N)$ are the pairs $(f,\phi)$
consisting of a morphism $f:(X,\underline M)\to(Y,\underline N)$
of $(S,\underline N)$-schemes, and a morphism $\phi:f^*\cG\to\cF$
of $\cO_{\!X}$-modules. With this notation, we claim that the rule :
\set\begin{equation}\label{eq_funct-Omega}
(X,\underline M)\mapsto
((X,\underline M),\Omega^1_{X/S}(\log\underline M/\underline N))
\end{equation}
defines a functor $\bprelog/(S,\underline N)\to\logMod/(S,\underline N)$.
Indeed, slightly more generally, consider a commutative diagram
of schemes with pre-log structures :
\set\begin{equation}\label{eq_will-be-cart}
{\diagram
(X,\underline M) \ar[r]^-g \ar[d]_f & (X',\underline M') \ar[d]^{f'} \\
(S,\underline N) \ar[r]^-h & (S',\underline N').
\enddiagram}
\end{equation}
An $\cO_{\!X}$-linear map :
\set\begin{equation}\label{eq_pul-bak-log-diff}
g^*\Omega^1_{X'/S'}(\log\underline M'/\underline N')\xrightarrow{dg}
\Omega_{X/S}(\log\underline M/\underline N)
\end{equation}
is the same as a natural transformation of functors :
\set\begin{equation}\label{eq_pull-back-der}
\cDer_{(S,\underline N)}((X,\underline M),\cF)\to
\cDer_{(S',\underline N')}((X',\underline M'),g_*\cF)
\end{equation}
on all $\cO_{\!X}$-modules $\cF$. The latter can be defined
as follows. Let $(\partial,\log\partial)$ be an $(S,\underline N)$-linear
derivation of $\underline M$ with values in $\cF$; then
we deduce morphisms :
$$
\partial':\cO_{\!X'}\xrightarrow{g^\natural}g_*\cO_{\!X}
\xrightarrow{g_*\partial}g_*\cF
\qquad
\log\partial':\underline M'\xrightarrow{\log g}g_*\underline M
\xrightarrow{g_*\log\partial}g_*\cF
$$
and it is easily seen that $(\partial',\log\partial')$
is a $(S',\underline N')$-linear derivation of $\underline M'$
with values in $g_*\cF$.

\sset\subsubsection{}\label{subsec_complex-of-log-di}
Consider two morphisms
$(X,\underline M)\xrightarrow{f}(Y,\underline N)
\xrightarrow{g}(Z,\underline P)$ of schemes with
pre-log structures. A direct inspection of the
definitions shows that :
$$
\cDer_{(Y,\underline N)}((X,\underline M),\cF)=
\Ker(\cDer_{(Z,\underline P)}((X,\underline M),\cF)\to
\cDer_{(Z,\underline P)}((Y,\underline N),f_*\cF))
$$
for every $\cO_{\!X}$-module $\cF$, whence a right
exact sequence of $\cO_{\!X}$-modules :
\set\begin{equation}\label{eq_complex-of-log-diff}
f^*\Omega^1_{Y/Z}(\log\underline N/\underline P)\xrightarrow{df}
\Omega^1_{X/Z}(\log\underline M/\underline P)\to
\Omega^1_{X/Y}(\log\underline M/\underline N)\to 0
\end{equation}
extending the standard right exact sequence for the usual
sheaves of relative differentials.

\begin{proposition}\label{prop_cart-diff}
Suppose that the diagram \eqref{eq_will-be-cart} is
cartesian. Then the map \eqref{eq_pul-bak-log-diff}
is an isomorphism.
\end{proposition}
\begin{proof} If \eqref{eq_will-be-cart} is cartesian,
$X$ is the scheme $X'\times_{S'}S$, and $\underline M$
is the push-out of the diagram:
$$
f^{-1}\underline N\leftarrow(f'\circ g)^{-1}\underline N'
\xrightarrow{\phi}g^*\underline M'.
$$
Suppose now that $\log\partial:\underline M'\to g_*\cF$
and $\partial:\cO_{\!X'}\to g_*\underline\cF$ define
a $f'$-linear derivation of $\underline M'$. By adjunction,
we deduce morphisms $\alpha:g^{-1}\cO_{\!X'}\to\cF$
and $\beta:g^{-1}\underline M'\to\cF$. By construction,
we have $\beta\circ\phi=0$, hence $\beta$ extends uniquely
to a morphism $\log\partial':\underline M\to\cF$ such that
$\log\partial'\circ\log f=0$. Likewise, $\alpha$ extends
by linearity to a unique $f$-linear derivation
$\partial:\cO_{\!X}\to\cF$.
One checks easily that $(\partial',\log\partial')$ is a
$f$-linear derivation of $\underline M$, and that every
$f$-linear derivation of $\underline M$ with values in $\cF$
is obtained in this fashion.
\end{proof}

\sset\subsubsection{}\label{subsec_left-adj}
The functor \eqref{eq_funct-Omega} admits a left adjoint.
Indeed, let $((X,\underline M\xrightarrow{\alpha}\cO_{\!X}),\cF)$
be any object of $\logMod/(S,\underline N)$; we define an
$(S,\underline N)$-scheme $(X\oplus\cF,\underline M\oplus\cF)$
as follows. $X\oplus\cF$ is the spectrum of the $\cO_{\!X}$-algebra
$\cO_{\!X}\oplus\cF$, whose multiplication law is given by
the rule :
$$
(s,t)\cdot(s',t'):=(ss',st'+s't)
\qquad
\text{for every local section $s$ of $\cO_{\!X}$ and $t$ of $\cF$}.
$$
Likewise, we define a composition law on the sheaf
$\underline M\oplus\cF$, by the rule :
$$
(m,t)\cdot(m',t'):=(mm',\alpha(m)\cdot t'+\alpha(m')\cdot t)
\qquad
\text{for every local section $m$ of $\underline M$ and $t$ of $\cF$}.
$$
Then $\underline M\oplus\cF$ is a sheaf of monoids, and
$\alpha$ extends to a pre-log structure
$\alpha\oplus\one_\cF:\underline M\oplus\cF\to\cO_{\!X}\oplus\cF$.
The natural map $\cO_{\!X}\to\cO_{\!X}\oplus\cF$ is a morphism
of algebras, whence a natural map of schemes
$\pi:X\oplus\cF\to X$, which extends to a morphism of
schemes with pre-log structures
$(X\oplus\cF,\underline M\oplus\cF)\to(X,\underline M)$,
by letting $\log\pi:\pi^*\underline M\to\underline M\oplus\cF$
be the map induced by the natural monomorphism (notice that
$\pi^*$ induces an equivalence of sites $X_\tau\isom(X\oplus\cF)_\tau$).

Now, let $(Y,\underline P)$ be any $(S,\underline N)$-scheme,
and :
$$
\phi:\cF:=((X,\underline M),\cF)\to
\Omega:=((Y,\underline P),\Omega^1_{Y/S}(\log\underline P/\underline N))
$$
a morphism in $\logMod/(S,\underline N)$. By definition, $\phi$
consists of a morphism $f:(X,\underline M)\to(Y,\underline P)$
and an $\cO_{\!X}$-linear map
$f^*\Omega^1_{Y/S}(\log\underline P/\underline N)\to\cF$,
which is the same as a $(S,\underline N)$-linear derivation :
$$
\partial:\cO_{\!Y}\to f_*\cF
\qquad
\log\partial:\underline P\to f_*\cF.
$$
In turns, the latter yields a morphism of $(S,\underline N)$-schemes :
\set\begin{equation}\label{eq_stupid-two}
(Y\oplus f_*\cF,\underline P\oplus f_*\cF)\to(Y,\underline P)
\end{equation}
determined by the map of algebras :
$$
\cO_{\!Y}\to\cO_{\!Y}\oplus f_*\cF
\qquad
s\mapsto(s,\partial s)\quad
\text{for every local section $s$ of $\cO_{\!Y}$}
$$
and the map of monoids :
$$
\underline P\mapsto\underline P\oplus f_*\cF
\qquad
p\mapsto(p,\log\partial p)\quad
\text{for every local section $p$ of $\underline P$}.
$$
Finally, we compose \eqref{eq_stupid-two} with the natural
morphism
$$
(X\oplus\cF,\underline M\oplus\cF)\to
(Y\oplus f_*\cF,\underline P\oplus f_*\cF)
$$
that extends $f$, to obtain a morphism
$g_\phi:(X\oplus\cF,\underline M\oplus\cF)\to(Y,\underline P)$.
We leave to the reader the verification that the rule
$\phi\mapsto g_\phi$ establishes a natural bijection :
$$
\Hom_{\logMod/(S,\underline N)}(\cF,\Omega)\isom
\Hom_{\bprelog/(S,\underline N)}
((X\oplus\cF,\underline M\oplus\cF),(Y,\underline P)).
$$

\sset\subsubsection{}
In the situation of definition \ref{def_log-Kahler}, let
$(\partial,\log\partial)$ be an $f$-linear derivation of
$\underline M$ with values in an $\cO_{\!X}$-module $\cF$.
Consider the map :
$$
\partial':\cO^\times_{\!X}\to\cF
\quad :\quad
u\mapsto u^{-1}\cdot\partial u
\quad
\text{for all local sections $u$ of $\cO^\times_{\!X}$}.
$$
By definition, $\partial'\circ\alpha:\alpha^{-1}\underline M\to\cF$
agrees with the restriction of $\log\partial$; in view of
the cocartesian diagram \eqref{eq_cokart}, we deduce that
$\log\partial$ extends uniquely to an $f$-linear derivation
$\log\partial^{\log}$ of $\underline M^{\log}$. Let
$f^{\log}:(X,\underline M^{\log})\to(Y,\underline N^{\log})$
be the map deduced from $f$; a similar direct verification
shows that $\log\partial^{\log}$ is a $f^{\log}$-linear
derivation. There follow natural identifications :
\set\begin{equation}\label{eq_nat-ident-pre-log}
\Omega^1_{X/Y}(\log(\underline{M}/\underline{N}))=
\Omega^1_{X/Y}(\log(\underline{M}^{\log}/\underline{N}))=
\Omega^1_{X/Y}(\log(\underline{M}^{\log}/\underline{N}^{\log})).
\end{equation}
Moreover, if $\underline{M}$ and $\underline{N}$ are log
structures, the natural map :
\set\begin{equation}\label{eq_nat-diff-onto}
\cO_{\!X}\otimes_\Z\log\underline{M}^\gp\to
\Omega^1_{X/Y}(\log(\underline{M}/\underline{N}))
\qquad a\otimes b\mapsto a\cdot d\log(b)
\end{equation}
is an epimorphism. Indeed, we have $da=d(a+1)$ for every
local section $a\in\cO_{\!X}(U)$ (for any {\'e}tale $X$-scheme
$U$), and locally on $X_\tau$, either $a$ or $1+a$ is
invertible in $\cO_{\!X}$ (this holds certainly on the stalks,
hence on appropriate small neighborhoods $U'\to U$);
hence $da$ lies in the image of \eqref{eq_nat-diff-onto}.

\begin{example}\label{ex_torique}
Let $R$ be a ring, $\phi:N\to M$ be any map of monoids,
and set :
$$
S:=\Spec\,R
\qquad
S[M]:=\Spec\,R[M]
\qquad
S[N]:=\Spec\,R[N].
$$
Also, let $f:\Spec(R,M)\to\Spec(R,N)$ be the morphism of log
schemes induced by $\phi$ (see \eqref{subsec_constant-log}).
With this notation, we claim that \eqref{eq_nat-diff-onto}
induces an isomorphism :
$$
\cO_{\!S[M]}\otimes_\Z\Coker\,\phi^\gp\isom
\Omega^1_{S[M]/S[N]}(\log M_{S[M]}^{\log}/N_{S[N]}^{\log})
$$
To see this, we may use \eqref{eq_complex-of-log-diff} to reduce
to the case where $N=\{1\}$. Next, notice that the functor
$$
\Mnd^o\to\bprelog/(S,\cO^\times_{\!S})
\quad :\quad M\mapsto\Spec(R,M)
$$
commutes with limits, and the same holds for the functor
\eqref{eq_funct-Omega}, since the latter is a right adjoint. Hence,
we may assume that $M$ is finitely generated; then lemma
\ref{lem_finite-pres}(i) further reduces to the case where
$M=\N^{\oplus n}$ for some integer $n\geq 0$, and even to the case
where $n=1$. Set $X:=S[\N]$; it is easy to see that a $S$-linear
derivation of $\N_X^{\log}$ with values in an
$\cO_{\!X}$-module $\cF$, is completely determined by a map of
additive monoids $\N\to\Gamma(X,\cF)$, and the latter is the same as
an $\cO_{\!X}$-linear map $\cO_{\!X}\to\cF$, whence the contention.
\end{example}

\begin{lemma}\label{lem_log-differ}
Let $f:(X,\underline{M})\to(Y,\underline{N})$ be a morphism
of schemes with quasi-coherent log structures.
Then :
\begin{enumerate}
\item
The $\cO_{\!X}$-module $\Omega^1_{X/Y}(\log\underline M/\underline N)$
is quasi-coherent.
\item
If $\underline M$ is coherent, $X$ is noetherian,
and $f:X\to Y$ is locally of finite type, then
$\Omega^1_{X/Y}(\log\underline M/\underline N)$ is a
coherent $\cO_{\!X}$-module.
\item
If both $\underline M$ and $\underline N$ are coherent,
and $f:X\to Y$ is locally of finite presentation, then
$\Omega^1_{X/Y}(\log\underline M/\underline N)$ is a
quasi-coherent $\cO_{\!X}$-module of finite presentation.
\end{enumerate}
\end{lemma}
\begin{proof} Applying the right exact sequence
\eqref{eq_complex-of-log-diff} to the sequence
$(X,\underline{M})\xrightarrow{f}(Y,\underline{N})\to
(Y,\cO^\times_{\!Y})$, we may easily reduce to
the case where $\underline N=\cO^\times_{\!Y}$ is the
trivial log structure on $Y$. In this case, $\underline N$
admits the chart given by the unique map of monoids :
$\{1\}\to\Gamma(Y,\cO_{\!Y})$, and $f$ admits the
chart $\{1\}\to P$, whenever $\underline M$ admits
a chart $P\to\Gamma(X,\cO_{\!X})$.

Hence, everything follows from the following
assertion, whose proof shall be left to the reader.
Suppose that $\phi:A\to B$ is a ring homomorphism,
$\underline M$ (resp. $\underline N$) is the constant
log structure on $X:=\Spec\,B$ (resp. on $Y:=\Spec\,A$)
associated to a map of monoids $\alpha:P\to B$ (resp.
$\beta:Q\to A$), and $f:(X,\underline M)\to(Y,\underline N)$
is defined by $\phi$ admits a chart $\vartheta:Q\to P$.
Then $\Omega^1_{X/Y}(\log\underline M/\underline N)$
is the quasi-coherent $\cO_{\!X}$-module $L^\sim$, associated to
the $B$-module $L:=(\Omega^1_{B/A}\oplus(B\otimes_\Z P^\gp))/R$,
where $R$ is the submodule generated by the elements
of the form $(0,1\otimes\log\vartheta(q))$ for all $q\in Q$,
and those of the form $(d\alpha(m),-\alpha(m)\otimes\log m)$,
for all $m\in M$.
\end{proof}

\sset\subsubsection{}\label{subsec_def-smooth}
Let us fix a log scheme $(Y,\underline N)$. To any pair of
$(Y,\underline N)$-schemes $X:=(X,\underline M)$,
$X':=(X',\underline M')$, we attach a contravariant functor
$$
\cH_Y(X',X):(X'_\tau)^o\to\Set
$$
by assigning, to every object $U$ of $X'_\tau$, the set of all
morphisms $(U,\underline M'_{|U})\to(X,\underline M)$ of
$(Y,\underline N)$-schemes. It is easily seen that $\cH_Y(X',X)$ is
a sheaf on $X'_\tau$. Any morphism $\phi:(X'',\underline
M'')\to(X',\underline M')$ (resp. $\psi:(X,\underline
M)\to(X'',\underline M)$) of $(Y,\underline N)$-schemes induces a
map of sheaves :
$$
\phi^*:\phi^*\cH_Y(X',X)\to\cH_Y(X'',X)\qquad \text{(resp.
$\psi_*:\cH_Y(X',X)\to\cH_Y(X',X'')$)}
$$
in the obvious way.

\begin{definition}\label{def_Exact-immersion}
With the notation of \eqref{subsec_def-smooth} :
\begin{enumerate}
\item
We say that a morphism $i:(T',\underline{L'})\to(T,\underline{L})$
of log schemes is a {\em closed immersion\/} (resp. an {\em exact
closed immersion}) if the underlying morphism of schemes $T'\to T$
is a closed immersion, and $\log
i:i^*\underline{L}\to\underline{L'}$ is an epimorphism (resp. an
isomorphism) of $T'_\tau$-monoids.
\item
We say that a morphism $i:(T',\underline{L'})\to(T,\underline{L})$
of log schemes is an {\em exact nilpotent immersion} if $i$
is an exact closed immersion, and the ideal
$\cI:=\Ker(\cO_{\!T}\to i_*\cO_{\!T'})$ is nilpotent.
\item
We say that a morphism $f:(X,\underline{M})\to(Y,\underline{N})$ of
log schemes is {\em formally smooth\/} (resp. {\em formally
unramified}, resp. {\em formally \'etale}) if, for every exact
nilpotent immersion $i:T'\to T$ of fine $(Y,\underline N)$-schemes,
the induced map of sheaves $i^*:i^*\cH_Y(T,X)\to\cH_Y(T',X)$ is an
epimorphism (resp. a monomorphism, resp. an isomorphism).
\item
We say that a morphism $f:(X,\underline{M})\to(Y,\underline{N})$ of
log schemes is {\em smooth\/} (resp. {\em unramified}, resp. {\em
\'etale}) if the underlying morphism $X\to Y$ is locally of finite
presentation, and $f$ is formally smooth (resp. formally unramified,
resp. formally \'etale).
\end{enumerate}
\end{definition}

\begin{example}\label{ex_graph-is-closed}
Let $(S,\underline P)$ be a log scheme, $(f,\log f):(X,\underline
M)\to(Y,\underline N)$ be a morphism of $S$-schemes, such that $Y$
is a separated $S$-scheme. The pair
$(\one_{(X,\underline M)},(f,\log f))$ induces a morphism
$$
\Gamma_{\!\!f}:(X,\underline M)\to(X',\underline M'):=
(X,\underline M)\times_S(Y,\underline N)
$$
the {\em graph\/} of $f$. Then it is easily seen that $\Gamma_{\!\!f}$
is a closed immersion of log schemes. Indeed, the morphism of schemes
underlying $\Gamma_{\!\!f}$ is a closed immersion (\cite[Ch.I,
5.4.3]{EGAI}) and it is easily seen that the morphism
$\log\Gamma_{\!\!f}:\Gamma_f^*\underline M'\to\underline M$ is an
epimorphism on the underlying sheaves of sets, so it is {\em a
fortiori\/} an epimorphism of $X_\tau$-monoids.
\end{example}

\begin{proposition}\label{prop_sorite-smooth}
Let $f:(X,\underline M)\to(Y,\underline N)$, $g:(Y,\underline
N)\to(Z,\underline P)$, and $h:(Y',\underline N')\to(Y,\underline
N)$ be morphisms of log schemes. Denote by $\bP$ either one of the
properties : "formally smooth", "formally unramified", "formally
\'etale", "smooth", "unramified", "\'etale". The following holds :
\begin{enumerate}
\item
If $f$ and $g$ enjoy the property $\bP$, then the same holds for
$g\circ f$.
\item
If $(f,\log f)$ enjoys the property $\bP$, then the same holds for
$$
(f,\log f)\times_{(Y,\underline N)}(Y',\underline N'):
(X,\underline M)\times_{(Y,\underline N)}(Y',\underline N')\to
(Y',\underline N').
$$
\item
Let $(j_\lambda:U_\lambda\to X~|~\lambda\in\Lambda)$ be a covering
family in $X_\tau$; endow $U_\lambda$ with the log structure
$(\underline M)_{|U_\lambda}$ and suppose that $f_\lambda:=(f\circ
j_\lambda,(\log f)_{|U_\lambda}): (U_\lambda,(\underline
M)_{|U_\lambda})\to(Y,\underline N)$ enjoys the property $\bP$, for
every $\lambda\in\Lambda$. Then $f$ enjoys the property $\bP$ as
well.
\item
An open immersion of log schemes is \'etale.
\item
A closed immersion of log schemes is formally unramified.
\end{enumerate}
\end{proposition}
\begin{proof}(i), (ii) and (iv) are left to the reader.

(iii): To begin with, if each $f_\lambda$ is locally of finite
presentation, the same holds for $f$, by \cite[Ch.IV, lemme
17.7.5]{EGA4}, hence we may assume that $\bP$ is either "formally
smooth" or "formally unramified". (Clearly, the case where $\bP$ is
"formally \'etale" will follow.)

Now, let $i:T'\to T$ be an exact nilpotent immersion of
$(Y,\underline N)$-schemes, and $\xi$ a $\tau$-point of $T'$. By
inspecting the definitions, it is easily seen that the stalk
$\cH_Y(T,X)_\xi$ is the union of the images of the stalks
$\cH_Y(T,U_\lambda)_\xi$, for every $\lambda\in\Lambda$, and
likewise for $\cH_Y(T',X)_\xi$. It readily follows that $f$ is
formally smooth whenever all the $f_\lambda$ are formally smooth.

Lastly, suppose that all the $f_\lambda$ are formally unramified,
and let $s_\xi,t_\xi\in\cH_Y(T,X)_{i(\xi)}$ be two sections whose
images in $\cH_Y(T',X)_\xi$ agree; after replacing $T$ by a
neighborhood of $i(\xi)$ in $T_\tau$, we may assume that
$s_\xi,t_\xi$ are represented by two $(Y,\underline N)$-morphisms
$s,t:T\to(X,\underline M)$ such that $s\circ i=t\circ i$. Choose
$\lambda\in\Lambda$ such that $U_\lambda$ is a neighborhood of
$s\circ i(\xi)=t\circ i(\xi)$; this means that there exist a
neighborhood $p':U'\to T'$ of $\xi$, and a morphism $s_{U'}:U'\to
U_\lambda$ lifiting $s\circ i$ (and thus, also $t\circ i$). Then we
may find a neighborhood $p:U\to T$ of $i(\xi)$ which identifies $p'$
with $p\times_T\one_{T'}$ (\cite[Ch.IV, Th.18.1.2]{EGA4}), and since
$j_\lambda$ is \'etale, we may furthermore find morphisms
$s_U,t_U:U\to U_\lambda$  such that $j_\lambda\circ s_U=s\circ
i=t\circ i=j_\lambda\circ t_U$. Set $i_U:=i\times_T\one_U:U'\to U$;
by construction, $s_U$ and $t_U$ yield two sections of
$i^*_U\cH_Y(U,U_\lambda)_\xi$, whose images in
$\cH_Y(U',U_\lambda)_\xi$ coincide. Since $f_\lambda$ is formally
unramified, it follows that -- up to replacing $U$ by a neighborhood
of $i(\xi)$ in $U_\tau$ -- we must have $s_U=t_U$, so $s_\xi=t_\xi$,
and we conclude that $f$ is formally unramified.

(v): Consider a commutative diagram of log schemes :
$$
\xymatrix{
(T',\underline L') \ar[r]^-{h'} \ar[d]_i & (X,\underline M) \ar[d]^f \\
(T,\underline L) \ar[r]^-h & (Y,\underline N) }
$$
where $f$ is a closed immersion, and $i$ is an exact closed
immersion. We are easily reduced to showing that there exists at
most a morphism $(g,\log ):(T,\underline L)\to(X,\underline M)$ such
that $f\circ g=h$ and $h'=g\circ i$. Since $f:X\to Y$ is a closed
immersion of schemes, there exists at most one morphism of schemes
$g:T\to X$ lifting $h$ and extending $h'$ (\cite[Ch.IV,
Prop.17.1.3(i)]{EGA4}). Hence we may assume that such a $g$ is
given, and we need to check that there exists at most one morphism
$\log g:g^*\underline M\to\underline L$ whose composition with
$g^*(\log f):h^*\underline N\to g^*\underline M$ equals $\log h$.
However, by assumption $\log f$ is an epimorphism, hence the same
holds for $g^*(\log f)$ (proposition \ref{prop_was-get-maddd}(iv)),
whence the contention.
\end{proof}

\begin{corollary}\label{cor_sorite-smooth}
Let $f$ and $g$ be as proposition {\em\ref{prop_sorite-smooth}}. We
have :
\begin{enumerate}
\item
If $g\circ f$ is formally unramified, the same holds for $f$.
\item
If $g\circ f$ is formally smooth (resp. formally \'etale) and $g$ is
formally unramified, then $f$ is formally smooth (resp. formally
\'etale).
\item
Suppose that $g$ is formally \'etale. Then $f$ is formally smooth
(resp. formally unramified, resp. formally \'etale) if and only if
the same holds for $g\circ f$.
\end{enumerate}
\end{corollary}
\begin{proof} (i): Let $Y=\bigcup_{\lambda\in\Lambda}$ be an affine
open covering of $Y$, and for each $\lambda\in\Lambda$, let
$f_\lambda:U_\lambda\times_Y(X,\underline M)\to(U_\lambda,\underline
N_{|U_\lambda})$ be the restriction of $f$; in light of proposition
\ref{prop_sorite-smooth}(iii) it suffices to show that each
$f_\lambda$ is formally unramified, and on the other hand, the
restriction $g\circ f_\lambda:U_\lambda\times_Y(X,\underline
M)\to(Z,\underline P)$ of $g\circ f$ is formally unramified, by
proposition \ref{prop_sorite-smooth}(i),(iv). It follows that we may
replace $f$ and $g$ respectively by $f_\lambda$ and $g_\lambda$,
which allows to assume that $Y$ is affine, especially separated, so
that $g$ is a separated morphism of schemes. In such situation, one
may -- in view of example \ref{ex_graph-is-closed} -- argue as in
the proof of \cite[Ch.IV, Prop.17.1.3]{EGA4} : the details shall be
left to the reader.

(ii) is a formal consequence of the definitions (cp. the proof of
\cite[Ch.IV, Prop.17.1.4]{EGA4}), and (iii) follows from (ii) and
proposition \ref{prop_sorite-smooth}(i).
\end{proof}

\begin{proposition}\label{prop_about-smooth-diff}
Let $f:(X,\underline M)\to(Y,\underline N)$ be a morphism of log
schemes, and $i$ an exact closed immersion of\/ $(Y,\underline
N)$-schemes, defined by an ideal $\cI:=\Ker(\cO_{\!T}\to
i_*\cO_{\!T'})$ with $\cI^2=0$. For any global section
$s:T'\to\cH_Y(T',X)$, denote by $\cT_s$ the morphism :
$$
i^*\cH_Y(T,X)\times_{\cH_Y(T',X)}T'\to T'
$$
deduced from $i^*:i^*\cH_Y(T,X)\to\cH_Y(T',X)$. Let also
$U\subset T'$ be the {\em image} of\/ $\cT_s$ ({\em i.e.} the
subset of all $t'\in T'$ such that $\cT_{s,\xi}\neq\emptyset$
for every $\tau$-point $\xi$ localized at $t'$), and suppose
that $U\neq\emptyset$. We have :
\begin{enumerate}
\item
$U$ is an open subset of\/ $T'$, and $\cT_{s|U}$ is a torsor
for the abelian sheaf
$$
\cG:=\cHom_{\cO_{\!T'}}
(s^*\Omega^1_{X/Y}(\log(\underline M/\underline N)),\cI)_{|U}.
$$
\item
If $f$ is a smooth morphism of log schemes with coherent
log structures, we have :
\begin{enumerate}
\item
The $\cO_{\!X}$-module $\Omega^1_{X/Y}(\log(\underline M/\underline N))$
is locally free of finite type.
\item
If\/ $T'$ is affine, $\cT_s$ is a trivial $\cG$-torsor.
\end{enumerate}
\end{enumerate}
\end{proposition}
\begin{proof}(i): To any $\tau$-point $\xi$ of $U$, and any two
given local sections $h$ and $g$ of $\cT_{s,\xi}$, we assign the
$f$-linear derivation of $\underline M_{s(\xi)}$ with values in
$s_*\cI_\xi$ given by the rule :
$$
\begin{aligned}
a & \mapsto h^*(a)-g^*(a) & &\text{for every $a\in\cO_{\!X,s(\xi)}$} \\
d\log m & \mapsto\log u(m):=u(m)-1 & &
\text{for every $m\in\underline M_{s(\xi)}$}
\end{aligned}
$$
where $u(m)$ is the unique local section of
$\Ker(\cO_{\!T,\xi}\to i_*\cO_{\!T',\xi})$ such that
$$
\log(g(m)\cdot u(m))=\log h(m).
$$
We leave to the reader the laborious, but straightforward
verification that the above map is well-defined, and yields
the sought bijection between $\cT_{s,\xi}$ and $\cG_\xi$.

(ii.b) is standard : first, since $f$ is smooth, we have $U=T'$; by
lemma \ref{lem_log-differ}(iii), the $\cO_{\!X}$-module
$\Omega^1_{X/Y}(\log(\underline M/\underline N))$ is quasi-coherent
and finitely presented, hence $\cG$ is quasi-coherent; however, the
obstruction to gluing local sections of a $\cG$-torsor lies in
$H^1(T_\tau,\cG)$ (see \eqref{subsec_coh_of_tors}); the latter
vanishes whenever $T$ is affine.

(ii.a) We may assume that $X$ is affine, say $X=\Spec\,A$;
then $\Omega^1_{X/Y}(\log(\underline M/\underline N))$ is the
quasi-coherent $\cO_{\!X}$-module arising from a finitely
presented $A$-module $\Omega$ (lemma \ref{lem_log-differ}(iii)).
Let $I$ be any $A$-module, and set :
$$
T_I:=X\oplus I^\sim
\qquad
\underline L_I:=\underline M\oplus I^\sim
$$
(notation of \eqref{subsec_left-adj}). Let $i:X\to T_I$ (resp.
$\pi:T_I\to X$) be the natural closed immersion (resp. the natural
projection); then $(T_I,\underline L_I)$ is a fine log scheme, and
$\pi$ (resp. $i$) extends to a morphism of log schemes
$(\pi,\log\pi):(T_I,\underline L_I)\to(X,\underline M)$ (resp.
$(i,\log i):(X,\underline M)\to(T_I,\underline L_I)$) with
$\log\pi:\pi^*\underline M\to\underline L_I$ (resp. $\log
i:i^*\underline L_I\to\underline M$) induced by the obvious
inclusion (resp. projection) map. Notice that $(i,\log i)$ is an
exact immersion and $(I^\sim)^2=0$, hence (i) says that the set
$\cT(T_I)$ of all morphisms $g:(T_I,\underline L_I)\to(X,\underline
M)$ such that
$$
f\circ g=f\circ\pi \qquad  g\circ i=\one_{(X,\underline M)}
$$
is in bijection with $\Hom_A(\Omega,I)$.
Moreover, any map $\phi:I\to J$ of $A$-modules induces a morphism
$i_\phi:(T_J,\underline L_J)\to(T_I,\underline L_I)$ of log schemes,
and if $\phi$ is surjective, $i_\phi$ is an exact nilpotent
immersion. Furthermore, we have a commutative diagram of sets :
$$
\xymatrix{ \cT(T_I) \ar[r]^-{i_\phi^*} \ar[d] & \cT(T_J) \ar[d] \\
           \Hom_A(\Omega,I) \ar[r]^-{\phi_*} & \Hom_A(\Omega,J)
}$$
whose vertical arrows are bijections, and where $i_\phi^*$
(resp. $\phi_*$) is given by the rule $g\mapsto g\circ i_\phi$,
(resp. $\psi\mapsto\phi\circ\psi$). Assertion (ii.b) implies
that $\iota_\phi^*$ is surjective when $f$ is smooth and
$\phi$ is surjective, hence the same holds for $\phi_*$,
{\em i.e.} $\Omega$ is a projective $A$-module, as stated.
\end{proof}

\begin{corollary}\label{cor_undercover}
Let $(f,\log f):(X,\underline M)\to(Y,\underline N)$
be a morphism of log schemes. We have:
\begin{enumerate}
\item
If $\underline M$ and $\underline N$ are fine log structures,
and $f$ is strict (see definition {\em\ref{def_amorph-log}(ii)}),
then $(f,\log f)$ is smooth (resp. {\'e}tale) if and only if
the underlying morphism of schemes $f:X\to Y$ is smooth
(resp. {\'e}tale).
\item
Suppose that $(f,\log f)$ is a smooth (resp. \'etale) morphism of
log schemes on the Zariski sites of $X$ and $Y$, and either :
\begin{enumerate}
\item
$\underline M$ and $\underline N$ are both integral log structures,
\item
or else $\underline N$ is coherent (on $Y_\Zar$) and $\underline M$
is fine (on $X_\Zar$).
\end{enumerate}
Then the induced morphism of log schemes on {\'e}tale sites :
$$
\tilde u{}^*(f,\log f):=(f,\tilde u{}^*_X\log f): \tilde
u{}^*_{\!X}(X,\underline M)\to\tilde u{}^*_Y(Y,\underline N)
$$
is smooth (resp. \'etale). (Notation of \eqref{subsec_choose-a-top}.)
\item
Suppose that $(f,\log f)$ is a morphism of log schemes on the
Zariski sites of $X$ and $Y$, and $\underline M$ is an integral log
structure on $X_\Zar$. Suppose also that $\tilde u{}^*(f,\log f)$ is
smooth (resp. {\'e}tale). Then the same holds for $(f,\log f)$.
\end{enumerate}
\end{corollary}
\begin{proof}(i): Suppose first that $(f,\log f)$ is smooth
(resp. {\'e}tale). Let $i:T'\to T$ be a nilpotent immersion of
affine schemes, defined by an ideal $\cI\subset\cO_T$ such that
$\cI^2=0$; let $s:T'\to X$ and $t:T\to Y$ be morphisms of schemes
such that $f\circ s=t\circ i$. By lemma
\ref{lem_simple-charts-top}(i), the log structures $\underline
L:=t^*\underline N$ and $\underline L':=s^*\underline M$ are fine,
and by choosing the obvious maps $\log s$ and $\log t$, we deduce a
commutative diagram :
\set\begin{equation}\label{eq_def-smooth}
{\diagram
(T',\underline{L'}) \ar[r]^s \ar[d]_{i} & (X,\underline{M}) \ar[d]^f \\
(T,\underline{L}) \ar[r]^t & (Y,\underline{N}).
\enddiagram}
\end{equation}
Then proposition \ref{prop_about-smooth-diff}(ii.b) says that there
exists a morphism (resp. a unique morphism) of schemes $g:T\to X$
such that $g\circ i=s$ and $f\circ g=t$, {\em i.e.} $f$ is smooth
(resp. \'etale).

The converse is easy, and shall be left as an exercise
for the reader.

(ii): Suppose first that $(f,\log f)$ is smooth, $\underline N$ is
coherent and $\underline M$ is fine. By proposition
\ref{prop_sorite-smooth}(iii), the assertion to prove is local on
$X_\et$, hence we may assume that $f$ admits a chart, given by a
morphism of finite monoids $P\to P'$ and commutative diagrams :
$$
P'_{X_\Zar}\to\underline M \qquad \omega:P_{Y_\Zar}\to\underline N
$$
(theorem \ref{th_good-charts}(i)). Now, consider a commutative
diagram of log schemes on {\'e}tale sites :
\set\begin{equation}\label{eq_etal-smooth}
{\diagram
(T',\underline{L'}) \ar[r]^-s \ar[d]_{i} &
\tilde u{}^*_X(X,\underline{M}) \ar[d]^{\tilde u{}^*f} \\
(T,\underline{L}) \ar[r]^-t & \tilde u{}^*_Y(Y,\underline{N})
\enddiagram}
\end{equation}
where $i$ is an exact nilpotent immersion of fine log schemes,
defined by an ideal $\cI\subset\cO_T$ such that $\cI^2=0$. Since the
assertion to prove is local on $T_\et$, we may assume that $T$ is
affine and -- in view of theorem \ref{th_good-charts}(ii) -- that
$t$ admits a chart, given by a morphism of finite monoids $\phi:P\to
Q$, and commutative diagrams :
\set\begin{equation}\label{eq_to_inspect} {\diagram
P_{T_\et}=t^*P_{Y_\et} \ar[rr]^-{t^*u{}^*_Y\omega} \ar[d]_{\phi_T} &
& t^*\tilde u{}^*_Y\underline N \ar[d]^{\log t} &
P \ar[r] \ar[d]_\phi & \Gamma(Y,\cO_Y) \ar[d]^{t^\natural} \\
Q_{T_\et} \ar[rr]^-\beta & & \underline L & Q \ar[r]^-\psi &
\Gamma(T,\cO_T) \enddiagram}
\end{equation}
Since $i$ is an exact closed immersion, it follows that the morphism
:
$$
Q_{T'_\et}\xrightarrow{ i^*\beta }i^*\underline L\xrightarrow{ \log
i }\underline L'
$$
is a chart for $\underline L'$. Especially, $(T',\underline L')$ is
isomorphic to $(T'_\et,Q^{\log}_{T'})$, the constant log structure
deduced from the morphism
$i^\natural\circ\psi:Q\to\Gamma(T',\cO_{T'})$. The latter is also
the log scheme $\tilde u{}^*_{T'}(T_\Zar,Q^{\log}_{T'})$. From
proposition \ref{prop_reduce-to-etale}(ii) we deduce that there
exists a unique morphism
$s_\Zar:(T_\Zar,Q^{\log}_{T'})\to(X,\underline M)$, such that
$\tilde u^*s_\Zar=s$. On the other hand, by inspecting
\eqref{eq_to_inspect} we find that there exists a unique morphism
$t_\Zar:(T_\Zar,Q^{\log}_T)\to(Y,\underline N)$ such that $\tilde
u{}^*t_\Zar=t$. These morphisms can be assembled into a diagram :
\set\begin{equation}\label{eq_Zariskified} {\diagram
(T'_\Zar,Q^{\log}_{T'})
\ar[d]_{i_\Zar} \ar[r]^-{s_\Zar} & (X,\underline M) \ar[d]^f \\
(T_\Zar,Q^{\log}_T) \ar[r]^-{t_\Zar} & (Y,\underline N) \enddiagram}
\end{equation}
where $i_\Zar$ is an exact closed immersion. By construction, the
diagram $\tilde u^*\eqref{eq_Zariskified}$ is naturally isomorphic
to \eqref{eq_etal-smooth}; especially, \eqref{eq_Zariskified}
commutes (proposition \ref{prop_reduce-to-etale}(i)). Now, since $f$
is smooth, proposition \ref{prop_about-smooth-diff}(ii.b) implies
that there exists a morphism
$v:(T_\Zar,Q^{\log}_T)\to(X,\underline M)$ such that $f\circ
v=t_\Zar$ and $v\circ i_\Zar=s_\Zar$; then $\tilde u{}^*v$ provides
an appropriate lifting of $t$, which allows to conclude that $\tilde
u{}^*(f,\log f)$ is smooth.

Next, suppose that $(f,\log f)$ is \'etale (and we are still in case
(b) of the corollary); then there exists a unique morphism $v$ with
the properties stated above. However, from proposition
\ref{prop_reduce-to-etale}(i),(ii) we deduce easily that the natural
map :
$$
\cH_Y((T_\Zar,Q^{\log}_T),(X,\underline M))(T)\to\cH_{\tilde
u{}^*(Y,\underline N)}((T,\underline L),\tilde u{}^*(X,\underline
M))(T)
$$
is a bijection, and the same holds for the analogous map for $T'$.
In view of the foregoing, this shows that the map $\cH_{\tilde
u{}^*(Y,\underline N)}((T,\underline L),\tilde u{}^*(X,\underline
M))(T)\to\cH_{\tilde u{}^*(Y,\underline N)}((T',\underline
L'),\tilde u{}^*(X,\underline M))(T')$ is bijective, whenever $t$
admits a chart and $T$ is affine. We easily conclude that $\tilde
u{}^*(f,\log f)$ is \'etale.

The case where both $\underline N$ and $\underline M$ are both
integral, is similar, though easier : we may assume that $\underline
L$ admits a chart, in which case $(T,\underline L)$ is of the form
$\tilde u{}^*(T_\Zar,Q^{\log}_T)$ for some finite integral monoid
$Q$; then $(T',\underline L')$ admits a similar description, and
again, by appealing to proposition \ref{prop_reduce-to-etale}(ii) we
deduce that \eqref{eq_etal-smooth} is of the form $\tilde
u^*\eqref{eq_Zariskified}$, in which case we conclude as in the
foregoing.

Conversely, if $\tilde u^*(f,\log f)$ is smooth, consider again a
commutative diagram \eqref{eq_def-smooth} of log schemes on
Zariski sites, with $i$ an exact closed immersion, defined by a
sheaf of ideals $\cI\subset\cO_{\!T}$ such that $\cI^2=0$. By
applying everywhere the pull-back functors from Zariski to \'etale
sites, we deduce a commutative diagram $\tilde
u^*\eqref{eq_def-smooth}$ of log schemes on \'etale sites, and it
is easy to see that $\tilde u{}^*(i,\log i)$ is again an exact
nilpotent immersion. According to proposition
\ref{prop_about-smooth-diff}(ii.b), after replacing $T$ by any
affine open subset, we may find a morphism of log schemes $g:\tilde
u{}^*_T(T,\underline L)\to \tilde u{}^*_X(X,\underline M)$, such
that $g\circ\tilde u{}^*i=\tilde u{}^*s$ and $\tilde u{}^* f\circ
g=t$. By proposition \ref{prop_reduce-to-etale}(i),(ii) there exists
a unique morphism $g':(T,\underline L)\to(X,\underline M)$ such that
$\tilde u{}^*g'=g$, and necessarily $g'\circ i=s$ and $f\circ g'=t$.
We conclude that $(f,\log f)$ is smooth.

Finally, if $\tilde u{}^*(f,\log f)$ is \'etale, the morphism $g$
exhibited above is unique, and therefore the same holds for $g'$, so
$(f,\log f)$ is \'etale as well.
\end{proof}

\begin{proposition}\label{prop_up-and-daun}
In the situation of \eqref{subsec_complex-of-log-di}, suppose
that the log structures $\underline M,\underline N,\underline P$
are coherent, and consider the following conditions :
\begin{enumerate}
\alphaenu
\item
$f$ is smooth (resp. {\'e}tale).
\item
$df$ is a locally split monomorphism (resp. an isomorphism).
\end{enumerate}
Then {\em (a)}$\Rightarrow${\em (b)}, and if $g\circ f$ is
smooth, then {\em (b)}$\Rightarrow${\em (a)}.
\end{proposition}
\begin{proof} We may assume that the schemes under consideration are
affine, say $X=\Spec\,A$, $Y=\Spec\,B$, $Z=\Spec\,C$; then
\eqref{eq_complex-of-log-diff} amounts to an exact sequence
of $A$-modules :
$$
A\otimes_B\Omega(g)\xrightarrow{df}\Omega(g\circ f)
\xrightarrow{\omega}\Omega(f)\to 0.
$$
On the other hand, let $i:(T'\underline L')\to(T,\underline L)$ be
an exact closed immersion of $(Y,\underline N)$-schemes, defined by
an ideal $\cI\subset\cO_{\!T}$ with $\cI^2=0$. Suppose that
$s:T'\to\cH_Y(T',X)$ is a global section, {\em i.e.} a given
morphism of $(Y,\underline N)$-schemes $(T'\underline
L')\to(X,\underline M)$. In this situation, we have a natural
sequence of morphisms :
$$
\cH_1:=\cH_Y(T,X)\xrightarrow{ \alpha
}\cH_2:=\cH_Z(T,X)\xrightarrow{ f_* }\cH_3:=\cH_Z(T,Y)
$$
where $\alpha$ is the obvious monomorphism. Let us consider the
pull-back of these sheaves along the global section $s$ :
$$
\cT:=T'\times_{\cH_Y(T',X)}i^*\cH_1 \qquad
\cT':=T'\times_{\cH_Z(T',X)}i^*\cH_2 \qquad
\cT'':=T'\times_{\cH_Z(T',Y)}i^*\cH_3.
$$
By proposition \ref{prop_about-smooth-diff}(i), any choice of a
global section of $\cT(T)$ determines a commutative diagram of
sets :
\set\begin{equation}\label{eq_haiving} {\diagram
\Hom_A(\Omega(f),\cI(T)) \ar[r]^-{\omega^*} \ar[d] &
\Hom_A(\Omega(g\circ f),\cI(T)) \ar[r]^-{df^*} \ar[d] &
\Hom_B(\Omega(g),\cI(T)) \ar[d] \\
\cT(T) \ar[r] & \cT'(T) \ar[r]^-{f_*} & \cT''(T)
\enddiagram}\end{equation}
whose vertical arrows are bijections.

In order to show that $df$ is a locally split monomorphism (resp. an
isomorphism), it suffices to show that the map $df^*:=\Hom_A(df,I)$
is surjective for every $A$-module $I$. To this aim, we take
$(T,L):=(T_I,\underline L_I)$, $(T',\underline L'):=(X,\underline
M)$, and let $i$ be the nilpotent immersion defined by the ideal
$I\subset A\oplus I$, as in the proof of proposition
\ref{prop_about-smooth-diff}(ii.a). If $\cI\subset\cO_{\!T}$ is the
corresponding sheaf of ideals, then $\cI(T)=I$; moreover, the
inclusions $A\to A\oplus I$ and $\underline M\subset\underline
M\oplus\cI$ determine a morphism $h:(T_I,\underline
L_I)\to(X,\underline M)$, which is a global section of $\cT(T)$.
Having made these choices, consider the resulting diagram
\eqref{eq_haiving} : if $f$ is {\'e}tale, $f_*$ is an isomorphism,
and when $f$ is smooth, $f_*$ is surjective (proposition
\ref{prop_about-smooth-diff}(ii.b)), whence the contention.

For the converse, we may suppose that $df$ is a split
monomorphism, hence $df^*$ in \eqref{eq_haiving} is a
split surjection. We have to show that $\cT(T)$ is not
empty, and by assumption (and proposition
\ref{prop_about-smooth-diff}(ii.b)) we know that
$\cT'(T)\neq\emptyset$. Choose any $\tilde h\in\cT'(T)$;
then $t$ and $f_*\tilde h$ are two elements of $\cT''(T)$,
so we may find $\phi\in\Hom_B(\Omega(g),\cI(T))$ such that
$\phi+f_*\tilde h=t$ (where the sum denotes the action
of $\Hom_B(\Omega(g),\cI(T))$ on its torsor $\cT''$).
Then we may write $\phi=\psi\circ df$ for some
$\psi:\Omega(g\circ f)\to\cI(Y)$, and it follows
easily that $\psi+\tilde h$ lies in $\cT(T)$.
Finally, if $df$ is an isomorphism, we have $\Omega(f)=0$,
hence $\cT(T)$ contains exactly one element.
\end{proof}

\begin{proposition}\label{prop_toric-smooth}
Let $R$ be a ring, $\phi:P\to Q$ a morphism of finitely generated
monoids, such that $\Ker\,\phi^\gp$ and the torsion subgroup
of\/ $\Coker\,\phi^\gp$ (resp. $\Ker\,\phi^\gp$ and
$\Coker\,\phi^\gp$) are finite groups whose orders are
invertible in $R$. Then, the induced morphism
$$
\Spec(R,\phi):\Spec(R,Q)\to\Spec(R,P)
$$
 is smooth (resp. \'etale).
\end{proposition}
\begin{proof} To ease notation, set
$$
f:=\Spec(R,\phi)
\qquad
(X,\underline M):=\Spec(R,Q)
\qquad
(Y,\underline N):=\Spec(R,P).
$$
Clearly $f$ is finitely presented. We have to show that, for
every commutative diagram like \eqref{eq_def-smooth} with $i$
an exact nilpotent immersion of fine log schemes, there is,
locally on $T_\tau$ at least one morphism (resp. a unique morphism)
$h:(T,\underline L)\to(X,\underline M)$ such that $f\circ h=h\circ i$.

Let $\cI:=\Ker(\cO_{\!T}\to i_*\cO_{\!T'})$; by considering
the $\cI$-adic filtration on $\cO_{\!T}$, we reduce easily
to the case where $\cI^2=0$, and then we may embed $\cI$ in
$\underline L$ via the morphism :
$$
\cI\to\cO^\times_{\!T}\subset\underline L\qquad x\mapsto 1+x.
$$
(Here $\cI$ is regarded as a sheaf of abelian groups, via
its addition law.) Since $i$ is exact, the natural morphism
$\underline L/\!\cI\to i_*\underline L'$ is an isomorphism,
whence a commutative diagram :
\set\begin{equation}\label{eq_would-Kato}
{\diagram
\underline L \ar[rr]^-{\log i} \ar[d] & & i_*\underline L' \ar[d] \\
\underline L^\gp \ar[rr]^-{(\log i)^\gp} & &
i_*(\underline L')^\gp\simeq\underline L^\gp/\cI
\enddiagram}
\end{equation}
and since $\underline L$ is integral, one sees easily that
\eqref{eq_would-Kato} is cartesian (it suffices to consider
the stalks over the $\tau$-points).

$\bullet$\ \
On the other hand, suppose first that both $\Ker\,\phi^\gp$
and $\Coker\,\phi^\gp$ are finite groups whose order is
invertible in $R$, hence in $\cI$; then a standard diagram
chase shows that we may find a unique map $g:P^\gp_T\to L^\gp$
of abelian sheaves that fits into a commutative diagram :
\set\begin{equation}\label{eq_many-names}
{\diagram
Q^\gp_T \ar[r] \ar[d]_{\phi_T^\gp} &
t^*\underline N^\gp \ar[rr]^-{(\log t)^\gp} & &
\underline L^\gp \ar[d]^{(\log i)^\gp} \\
P^\gp_T \ar[r] \ar[rrru]^g &
i_*s^*\underline M^\gp \ar[rr]^-{i_*(\log s)^\gp} & &
i_*(\underline L')^\gp.
\enddiagram}
\end{equation}

$\bullet$\ \ 
More generally, we may write : $\Coker\,\phi^\gp\simeq G\oplus H$,
where $H$ is a finite group with order invertible in $R$, and $G$
is a free abelian group of finite rank. The direct summand $G$
lifts to a direct summand $G'\subset P^\gp$. Extend $\phi^\gp$
to a map $\psi:Q^\gp\oplus G'\to P^\gp$, by the rule :
$(x,g)\mapsto\phi^\gp(x)\cdot g$. Given any $\tau$-point $\xi$
of $T$, we may extend the morphism
$Q^\gp=Q^\gp_{T,\xi}\to\underline L^\gp_\xi$ in
$\eqref{eq_many-names}_\xi$, to a map
$\omega:Q^\gp\oplus G'\to\underline L^\gp_\xi$ whose
composition with $(\log i)^\gp_\xi$ agrees with the
composition of $\psi$ and the bottom map
$P^\gp=P^\gp_{T,\xi}\to i_*(\underline L')^\gp_\xi$ of
$\eqref{eq_many-names}_\xi$. By the usual arguments,
$\omega$ extends to a map of abelian sheaves
$\vartheta:(Q^\gp\oplus G')_U\to(\underline L^\gp)_{|U}$
on some neighborhood $U\to T$ of $\xi$, and if $U$ is small
enough, the composition $(\log i)^\gp_{|U}\circ\vartheta$
agrees with the composition of $\psi_U$ and the bottom
map $\beta:P^\gp_U\to i_*(\underline L'^\gp)_{|U}$ of
$\eqref{eq_many-names}_{|U}$. We may then replace $T$
by $U$, and since $\Ker\,\psi=\Ker\,\phi^\gp$,
and $\Coker\,\psi=H$, the same diagram chase as in the
foregoing shows that we may again find a morphism
$g:P^\gp_T\to\underline L^\gp$ fitting into a commutative
diagram :
$$
\xymatrix{ (Q^\gp\oplus G')_T \ar[r]^-\vartheta \ar[d]_{\psi_T} &
\underline L^\gp \ar[d]^{(\log i)^\gp} \\
P^\gp_T \ar[r]^-\beta \ar[ru]^-g & i_*\underline L'^\gp
}$$
In either case, in view of \eqref{eq_would-Kato}, the morphisms
$$
P_T\to P^\gp_T\xrightarrow{g}\underline L^\gp
\qquad\text{and}\qquad
P_T\to i_*s^*\underline M\xrightarrow{i_*\log s}i_*\underline L'
$$
determine a unique morphism $P_T\to\underline L$, which
induces a morphism of log schemes
$(T,\underline L)\to(X,\underline M)$ with the sought property.
\end{proof}

\begin{theorem}\label{th_charact-smoothness}
Let $f:(X,\underline M)\to(Y,\underline N)$ be a morphism of fine
log schemes. Assume we are given a fine chart
$\beta:Q_Y\to\underline N$ of $\underline N$. Then the following
conditions are equivalent :
\begin{enumerate}
\alphaenu
\item
$f$ is smooth (resp. \'etale).
\item
There exist a covering family $(g_\lambda:U_\lambda\to
X~|~\lambda\in\Lambda)$ in $X_\et$, and for every
$\lambda\in\Lambda$, a fine chart $(\beta,(P_\lambda)_{U_\lambda}\to
g^*_\lambda\underline M,\phi_\lambda:Q\to P_\lambda)$ of the induced
morphism of log schemes
$$
f_{|U_\lambda}:=(f\circ g_\lambda,g^*_\lambda\log f):
(U_\lambda,g^*_\lambda\underline M)\to(Y,\underline N)
$$
such that :
\begin{enumerate}
\romanenuii
\item
$\Ker\,\phi^\gp_\lambda$ and the torsion subgroup of
$\Coker\,\phi_\lambda^\gp$ (resp. $\Ker\,\phi^\gp_\lambda$ and
$\Coker\,\phi^\gp_\lambda$) are finite groups of orders invertible
in $\cO_{\!U_\lambda}$.
\item
The natural morphism of\/ $Y$-schemes $p_\lambda:U_\lambda\to
Y\times_{\Spec\,\Z[Q]}\Spec\,\Z[P_\lambda]$ is {\'e}tale.
\end{enumerate}
\end{enumerate}
\end{theorem}
\begin{proof} Suppose first that $\tau=\et$, so $\underline M$
and $\underline N$ are log structures on {\'e}tale sites. Then the
log structure $g^*_\lambda\underline M$ on $(U_\lambda)_\et$ is just
the restriction of $\underline M$, and $g^*_\lambda\log f$ is the
restriction of $\log f$ to $(U_\lambda)_\et$. In case $\tau=\Zar$,
the log structure $g^*_\lambda\underline M$ can be described as
follows. Form the log scheme $(X,\underline M'):=\tilde
u{}^*_X(X,\underline M)$ (notation of \eqref{subsec_choose-a-top}),
take the restriction $\underline M'{}_{|U_\lambda}$ of $\underline M'$
to $(U_\lambda)_\et$, and push forward to the Zariski site to obtain
$\tilde u_{X*}(\underline M'{}_{|U_\lambda})=g^*_\lambda\underline M$
(by proposition \ref{prop_reduce-to-etale}(ii)). Since $f$ is smooth
(resp. {\'e}tale) if and only if $\tilde u^*f$ is (corollary
\ref{cor_undercover}(ii),(iii)), we conclude that the assertion
concerning $f$ holds if and only if the corresponding assertion for
$\tilde u^*f$ does. Hence, it suffices to consider the case of log
structures on {\'e}tale sites. Set $S:=\Spec\,\Z[Q]$.

(b)$\Rightarrow$(a): Taking into account lemma \ref{lem_time}(i),
we deduce a commutative diagram of log schemes :
$$
\xymatrix{
(U_\lambda,\underline M_{|U_\lambda})
\ar[r]^-\sim \ar[d]_{g_\lambda} &
\Spec(\Z,P_\lambda)\times_{\Spec\,\Z[P_\lambda]}U_\lambda
\ar[r]^-{p_\lambda} &
\Spec(\Z,P_\lambda)\times_{\Spec\,\Z[Q]}Y \ar[d]^{\pi_\lambda} \\
(X,\underline M) \ar[r]^-f & (Y,\underline N) \ar[r]^-\sim
& \Spec(\Z,Q)\times_{\Spec\,\Z[Q]}Y.
}$$
It follows by corollary \ref{cor_undercover}(i) (resp. by
propositions \ref{prop_sorite-smooth}(ii) and
\ref{prop_toric-smooth}) that $p_\lambda$ (resp. $\pi_\lambda$) is
smooth. Hence $f\circ g_\lambda$ is smooth, by proposition
\ref{prop_sorite-smooth}(i). Finally, $f$ is smooth, by proposition
\ref{prop_sorite-smooth}(iii).

(a)$\Rightarrow$(b): Suppose that $f$ is smooth, and fix
a geometric point $\xi$ of $X$. Since
$\eqref{eq_nat-diff-onto}_\xi$ is a surjection, we may
find elements $t_1,\dots,t_r\in\underline M{}_\xi$
such that $(d\log t_i~|~i=1,\dots,r)$ is a basis of the
free $\cO_{\!X,\xi}$-module
$\Omega^1_{X/Y}(\log\underline M/\underline N)_\xi$
(proposition \ref{prop_about-smooth-diff}(ii.a)). Moreover,
the kernel of $\eqref{eq_nat-diff-onto}_\xi$ is generated
by sections of the form :
\begin{itemize}
\item
$1\otimes \log a$ where
$a\in N':=
\Img(\log f_\xi:\underline N{}_{f(\xi)}\to\underline M{}_\xi)$
\item
$\sum_{j=1}^s\alpha(m_i)\otimes\log m_i$ where
$m_1,\dots,m_s\in\underline M{}_\xi$ and $\sum^s_{j=1}d\alpha(m_j)=0$
in $\Omega^1_{X/Y}$
\end{itemize}
whence a well-defined $\cO_{\!X,\xi}$-linear map :
\set\begin{equation}\label{eq_surge}
\Omega^1_{X/Y}(\log\underline M/\underline N)_\xi\to
\kappa(\xi)\otimes_\Z
(\underline M{}^\gp_\xi/(\underline M{}^\times_\xi\cdot N'))
\quad :\quad
d\log a\mapsto 1\otimes a.
\end{equation}
Consider the map of monoids :
$$
\phi:P_1:=\N^{\oplus r}\oplus Q\to\underline M{}_\xi
$$
which is given by the rule : $e_i\mapsto t_i$ on the canonical
basis $e_1,\dots,e_r$ of $\N^{\oplus r}$, and on the summand
$Q$ it is given by the map
$Q\xrightarrow{\beta_\xi}\underline N{}_{f(\xi)}
\xrightarrow{\log f_\xi}\underline M{}_\xi$.
Since \eqref{eq_surge} is a surjection, we see that
the same holds for the induced map
$$
\kappa(\xi)\otimes_\Z P_1^\gp
\xrightarrow{\one_{\kappa(\xi)}\otimes_\Z\phi^\gp}
\kappa(\xi)\otimes_\Z\underline M{}^\gp_\xi\to
\kappa(\xi)\otimes_\Z(\underline M{}^\gp_\xi/\underline M{}^\times_\xi).
$$
It follows that the cokernel of the map
$\bar\phi:P_1^\gp\to\underline M{}^\gp_\xi/\underline M{}^\times_\xi$
induced by $\phi^\gp$, is a finite group (lemma
\ref{lem_simple-charts}(i)) annihilated by an integer $n$ which is
invertible in $\cO_{\!X,\xi}$. Let
$m_1,\dots,m_s\in\underline M{}^\gp_\xi$ be finitely many elements,
whose images in $\underline M{}^\gp_\xi/\underline M{}^\times_\xi$
generate $\Coker\,\bar\phi$; therefore we may find
$u_1,\dots,u_s\in\underline M{}^\times_\xi$, and
$x_1,\dots,x_s\in P^\gp_1$, such that $m_i^n\cdot u_i=\bar\phi(x_i)$
for every $i\leq s$. However, since $\cO_{\!X,\xi}$ is strictly
henselian, $\underline M{}_\xi^\times\simeq\cO^\times_{\!X,\xi}$ is
$n$-divisible, hence we may find $v_1,\dots,v_s$ in
$\underline M{}_\xi^\times$ such that $u_i=v_i^n$ for $i=1,\dots,n$.
Define group homomorphisms :
$$
\Z^{\oplus s}\xrightarrow{\gamma}\Z^{\oplus s}\oplus P_1^\gp
\xrightarrow{\delta}\underline M{}^\gp_\xi
$$
by the rules : $\gamma(e_i)=(-ne_i,x_i)$ and
$\delta(e_i,y)=m_iv_i\cdot\phi(y)$ for every $i=1,\dots,s$ and every
$y\in P^\gp_1$. It is easily seen that $\delta$ factors through a
group homomorphism $h:G:=\Coker\,\gamma\to\underline M{}_\xi$;
moreover the natural map $P_1^\gp\to G$ is injective, and its
cokernel is annihilated by $n$; furthermore, the induced map
$G\to\underline M{}^\sharp_\xi$ is surjective. Let
$P:=h^{-1}\underline M{}_\xi$. Then, the natural map $Q^\gp\to P^\gp$
is injective, and the torsion subgroup of $P^\gp/Q^\gp$ is
annihilated by $n$. We deduce that the rule : $x\mapsto d\log h(x)$
for every $x\in P^\gp$, induces an isomorphism :
\set\begin{equation}\label{eq_stalk-ome-P}
\cO_{\!X,\xi}\otimes_\Z(P^\gp/Q^\gp)\isom
\Omega^1_{X/Y}(\log\underline M/\underline N)_\xi.
\end{equation}
It follows that we may find an {\'e}tale neighborhood
$U\to X$ of $\xi$, such that \eqref{eq_stalk-ome-P} extends
to an isomorphism of $\cO_{\!U}$-modules :
\set\begin{equation}\label{eq_extend-to-yu}
\cO_{\!U}\otimes_\Z(P^\gp/Q^\gp)\isom
\Omega^1_{U/Y}(\log\underline M{}_{|U}/\underline N).
\end{equation}
Next, proposition \ref{prop_good-charts} says that, after replacing
$U$ by a smaller \'etale neighborhood of $\xi$, the restriction
$h_{|P}:P\to\underline M{}_\xi$ extends to a chart
$P_U\to\underline M{}_{|U}$, whence a strict morphism
$$
p:(U,\underline M{}_{|U})\to(Y',P^{\log}_{Y'}):=
(Y,\underline N)\times_{\Spec(\Z,Q)}\Spec(\Z,P).
$$
as sought. Taking into account corollary \ref{cor_undercover}(i),
it remains only to show :
\begin{claim} $p$ is \'etale.
\end{claim}
\begin{pfclaim}[] By proposition \ref{prop_cart-diff} and example
\ref{ex_torique}, we have natural isomorphisms
$$
\Omega^1_{Y'/Y}(\log P^{\log}_{Y'}/\underline N)
\isom\cO_{\!Y'}\otimes_\Z(P^\gp/Q^\gp).
$$
In view of \eqref{eq_extend-to-yu}, it follows that
the map
$$
dp:p^*\Omega^1_{Y'/Y}(\log P^{\log}_{Y'}/\underline N)
\to\Omega^1_{U/Y}(\log\underline M{}_{|U}/\underline N)
$$
is an isomorphism, and then the claim follows from
proposition \ref{prop_up-and-daun}.
\end{pfclaim}
\end{proof}

\begin{corollary}\label{cor_charact-smoothness}
Keep the notation of theorem {\em\ref{th_charact-smoothness}}.
Suppose that $f$ is a smooth morphism of fs log schemes, and that
$Q$ is fine, sharp and saturated. Then there exist a covering
family $(g_\lambda:U_\lambda\to X)$ in $X_\et$, and fine and
saturated charts $(\beta,(P_\lambda)_{U_\lambda}\to
g^*_\lambda\underline M,\phi_\lambda:Q\to P_\lambda)$ of
$f_{|U_\lambda}$ fulfilling conditions {\em(b.i)} and {\em(b.ii)} of
the theorem, and such that moreover $\phi_\lambda$ is injective, and
$P^\times_\lambda$ is a torsion-free abelian group, for every
$\lambda\in\Lambda$.
\end{corollary}
\begin{proof} Notice that, under the stated assumptions, $Q^\gp$ is
a torsion-free abelian group; hence theorem
\ref{th_charact-smoothness} already implies the existence of a
covering family $(g_\lambda:U_\lambda\to X)$ in $X_\et$, and of fine
charts $(\beta,\omega'_\lambda:(P'_\lambda)_{U_\lambda}\to
g^*_\lambda\underline M,\phi'_\lambda:Q\to P'_\lambda)$ such that
$\phi'_\lambda$ is injective and $\Coker(\phi'_\lambda)^\gp$ is a
finite group of order invertible in $\cO_{U_\lambda}$, for every
$\lambda\in\Lambda$. Since $g_\lambda^*\underline M$ is saturated
(lemma \ref{lem_simple-charts-top}(i)), it is clear that
$\omega'_\lambda$ factors through a morphism
$(P'_\lambda)^\sat_{U_\lambda}\to g^*_\lambda\underline M$ of
pre-log structures, which is again a chart, so we may assume that
$P'_\lambda$ is fine and saturated, for every $\lambda\in\Lambda$
(corollary \ref{cor_fragment-Gordon}(ii)), in which case we may
find an isomorphism of monoids $P'_\lambda=P_\lambda\times G$, where
$P_\lambda^\times$ is torsion-free, and $G$ is a finite group (lemma
\ref{lem_decomp-sats}). Let $d$ be the order of
$\Coker(\phi'_\lambda)^\gp$, and denote by $\phi_\lambda$ the
composition of $\phi'_\lambda$  and the projection $P'_\lambda\to
P_\lambda$; since $Q^\gp$ is a torsion-free abelian group, we deduce
a short exact sequence :
\set\begin{equation}\label{eq_shorts}
0\to G\to\Coker(\phi'_\lambda)^\gp\to\Coker\,\phi_\lambda^\gp\to 0
\end{equation}
and notice that $\phi_\lambda$ is also injective. We may assume that
$U_\lambda$ and $Y$ are affine, say $U_\lambda=\Spec\,B_\lambda$ and
$Y=\Spec\,A$, and since $d$ is invertible in $\cO_{U_\lambda}$, we
reduce easily -- via base change by a finite morphism $Y'\to Y$ --
to the case where $A$ contains the subgroup
$\bmu_d\subset\bar\Q{}^\times$ of $d$-th power roots of $1$. The
chart $\omega'_\lambda$ determines a morphism of monoids
$P'_\lambda\to B_\lambda$, and the map $f^\sharp:A\to B_\lambda$
factors through the natural ring homomorphism
$$
A\to A\otimes_{R[Q]}R[P'_\lambda]\isom
A\otimes_{R[Q]}(R[P_\lambda]\otimes_R R[G]) \qquad\text{where}\qquad
R:=\Z[d^{-1},\bmu_d].
$$
On the other hand, let $\Gamma:=\Hom_\Z(G,\bmu_d)$; then we have a
natural decomposition
$$
R[G]\simeq\prod_{\chi\in\Gamma}e_\chi R[G]
$$
where $e_\chi$ is the idempotent of $R[G]$ defined as in
\eqref{eq_central-idemp} (cp. the proof of theorem \ref{th_horiz}(i));
each factor is a ring isomorphic to $R$, whence a corresponding
decomposition of $U_\lambda$ as a disjoint union of $Y$-schemes
$U_\lambda=\coprod_{\chi\in\Gamma}U_{\lambda,\chi}$. We are then
further reduced to the case where $U_\lambda=U_{\lambda,\chi}$ for
some character $\chi$ of $G$. In view of \eqref{eq_shorts}, it is
easily seen that the composition
$$
Q^\gp\xrightarrow{ (\phi'_\lambda)^\gp }(P'_\lambda)^\gp\to
G\xrightarrow{ \chi }\bmu_d
$$
extends to a well defined group homomorphism
$\bar\chi:P_\lambda^\gp\to\bmu_d$, whence a map
$\bar\chi_{U_\lambda}:P_{\lambda,U_\lambda}\to\mu_{d,U_\lambda}\subset
g^*_\lambda\underline M$ of sheaves on $U_{\lambda,\tau}$. Define
$\omega_\lambda:P_{\lambda,U_\lambda}\to g_\lambda^*\underline M$ by
the rule : $s\mapsto\omega'_\lambda(s)\cdot\bar\chi_{U_\lambda}(s)$
for every local section $s$ of $P_{\lambda,U_\lambda}$. It is easily
seen that $\omega_\lambda$ is again a chart for
$g^*_\lambda\underline M$ ({\em e.g.} one may apply lemma
\ref{lem_check-iso}). Lastly, a direct inspection shows that
$(\beta,\omega_\lambda,\phi_\lambda)$ is a chart of $f_{|U_\lambda}$
with the sought properties.
\end{proof}

\sset\subsubsection{}\label{subsec_annoying-gen-II}
Let $(Y_i~|~i\in I)$ be a cofiltered system of quasi-compact
and quasi-separated schemes, with affine transition morphisms,
and suppose that $0$ is an initial object of the indexing
category $I$. Suppose also that $g_0:X_0\to Y_0$ is a finitely
presented morphism of schemes. Let $X_i:=X_0\times_{Y_0}Y_i$
for every $i\in I$, and denote $g_i:X_i\to Y_i$ the induced
morphism. Let also $g:X\to Y$ be the limit of the family of
morphisms $(g_i~|~i\in I)$.

\begin{corollary}\label{cor_smooth-charct}
In the situation of \eqref{subsec_annoying-gen-II},
suppose that $(g,\log g):(X,\underline M)\to(Y,\underline N)$
is a smooth morphism of fine log schemes. Then there exist
$i\in I$, and a smooth morphism
$(g_i,\log g_i):(X_i,\underline M{}_i)\to(Y_i,\underline N{}_i)$
of fine log schemes, such that $\log g=\pi_i^*\log g_i$.
\end{corollary}
\begin{proof} First, using corollary \ref{cor_notime-now},
we may find $i\in I$ such that $(g,\log g)$ descends to a
morphism $(g_i,\log g_i)$ of log schemes with coherent
log structures. After replacing $I$ by $I/i$, we
may then suppose that $i=0$, in which case we set
$(X_i,\underline M{}_i):=X_i\times(X_0,\underline M{}_0)$,
and define likewise $(Y_i,\underline N{}_i)$ for every $i\in I$.

Next, arguing as in the proof of corollary
\ref{cor_descend-chart-from-infty}(ii), we may assume that
$\underline N{}_0$ admits a fine chart.
In this case, also $\underline N$ admits a fine chart,
and then we my find a covering family
$\cU:=(U_\lambda\to X~|~\lambda\in\Lambda)$ for $X_\tau$, such
that the induced morphism
$(U_\lambda,\underline M{}_{|U_\lambda})\to(Y,\underline N)$
admits a chart fulfilling conditions (i) and (ii) of theorem
\ref{th_charact-smoothness}.
Moreover, under the current assumptions, $X$ is quasi-compact,
hence we may assume that $\Lambda$ is a finite set,
in which case there exists $i\in I$ such that
$\cU$ descends to a covering family
$(U_{i,\lambda}\to X_i~|~\lambda\in\Lambda)$ for $X_{i,\tau}$
(claim \ref{cl_combination}(ii)), and as usual, we may then
reduce to the case where $i=0$. It then suffices to show that
there exists $i\in I$ such that the induced morphism
$U_{0,\lambda}\times_{X_0}(X_i,\underline M{}_i)\to
(Y_i,\underline N{}_i)$ is smooth (proposition
\ref{prop_sorite-smooth}(iii)). Thus, we may replace the system
$(X_i~|~i\in I)$ by the system of schemes
$(U_{0,\lambda}\times_{X_0}X_i~|~i\in I)$, and assume from start
that $(g,\log g)$ admits a fine chart fulfilling conditions (i)
and (ii) of theorem \ref{th_charact-smoothness}. In this case,
corollary \ref{cor_descend-chart-from-infty}(i) and
\cite[Ch.IV, Prop.17.7.8(ii)]{EGAIV-3} imply more precisely
that there exists $i\in I$ such that the morphism
$(X_i,\underline M{}_i)\to(Y_i,\underline N{}_i)$ fulfills
conditions (i) and (ii) of theorem \ref{th_charact-smoothness},
whence the contention.
\end{proof}

\subsection{Logarithmic blow up of a coherent ideal}
This section introduces the logarithmic version of
the scheme-theoretic blow up of a coherent ideal,
which will be exhibited as the logarithmic homogeneous
spectrum of a certain graded algebra, naturally attached
to any sheaf of ideals in a log structure.

\sset\subsubsection{}
We shall consider first the logarithmic counterparts of
the notions introduced in section \ref{sec_coh-blow-up}.
To begin with, let $X$ be any scheme, $N$ a monoid, and
$\underline P$ a (commutative) $N$-graded monoid of the
topos $X_\tau$ (see definition \ref{def_grad-monoids}).
Notice that $\underline P^\times=
\coprod_{n\in N}(\underline P^\times\cap\underline P_n)$,
hence the sheaf of invertible sections of a $N$-graded
monoid on $X$ is a $N$-graded abelian sheaf.

\begin{definition}\label{def_qcoh-log-gr}
Let $X$ be a scheme, and $\beta:\underline M\to\cO_X$ a log
structure on $X_\tau$.
\begin{enumerate}
\item
A {\em $\N$-graded $\cO_{\!(X,\underline M)}$-algebra} is a
datum $(\cA,\underline P,\alpha,\beta_{\!\cA})$ consisting of a
$\N$-graded $\cO_{\!X}$-algebra $\cA:=\oplus_{n\in\N}\cA_n$ (on
the site $X_\tau$), a graded monoid $\underline P$ on $X_\tau$,
and a commutative diagram :
$$
\xymatrix{ \underline M \ar[r]^-\alpha \ar[d]_\beta &
\underline P \ar[d]^{\beta_{\!\cA}} \\
\cO_{\!X} \ar[r] & \cA
}$$
where $\beta_{\!\cA}$ restricts to a morphism of
graded monoids $\underline P\to\coprod_{n\in\N}\cA_n$ (and
the composition law on the target is induced by the multiplication
law of $\cA$), $\alpha$ is a morphism of monoids
$\underline M\to(\underline P)_0$, and the bottom map is the
natural morphism $\cO_{\!X}\to\cA_0$. We say that the $\N$-graded
$\cO_{\!(X,\underline M)}$-algebra
$(\cA,\underline P,\alpha,\beta_{\!\cA})$ is {\em quasi-coherent},
if $\cA$ is a quasi-coherent $\cO_{\!X}$-algebra.
\item
A {\em morphism
$(g,\log g):(\cA,\underline P,\alpha,\beta_{\!\cA})\to
(\cA',\underline P',\alpha',\beta'_{\!\cA})$ of\/ $\N$-graded
$\cO_{\!(X,\underline M)}$-alge\-bras\/} is a commutative diagram :
$$
\xymatrix{ \underline P \ar[d]_{\beta_{\!\cA}} \ar[rr]^-{\log g}
& & \underline P' \ar[d]^{\beta'_{\!\cA}} \\
\cA \ar[rr]^-g & & \cA'
}$$
where $\log g$ is a morphism of $\N$-graded monoid such
that $\log g\circ\alpha=\alpha'$, and $g$ is a morphism
of $\cO_{\!Z}$-algebras.
\end{enumerate}
\end{definition}

\sset\subsubsection{}\label{subsec_inv-mod-mon}
Let $M$ be a monoid, $S$ an $M$-module; we say that an element
$s\in S$ is {\em invertible}, if the translation map
$M\to S$ : $m\mapsto m\cdot s$ (for all $m\in M$) is an
isomorphism. It is easily seen that the subset $S^\times$ consisting
of all invertible elements of $S$, is naturally an $M^\times$-module.

Let $\underline M$ be a sheaf of monoids on $X_\tau$, and $\cN$ an
$\underline M$-module; by restriction of scalars, $\cN$ is naturally
an $\underline M^\times$-module. We define the $\underline
M^\times$-submodule $\cN^\times\subset\cN$, by the rule :
$$
\cN^\times(U):=\cN(U)^\times
\qquad
\text{for every object $U$ of $X_\tau$}.
$$
Conversely, for any $\underline M^\times$-module $\cP$, the
extension of scalars $\cN\otimes_{\underline M^\times}\underline M$
defines a functor $\underline M^\times\Mod\to\underline M\Mod$. It
is easily seen that the latter restricts to an equivalence from the
full subcategory $\underline M^\times\text{-}\mathbf{Inv}$ of
$\underline M^\times$-torsors to the subcategory $(\underline
M\text{-}\mathbf{Inv})^\times$ of $\underline M\Mod$ whose objects
are all invertible $\underline M$-modules, and whose morphisms are
the isomorphisms of $\underline M$-modules (see definition
\ref{def_coh-idea-log}(iv)); the functor $\cN\mapsto\cN^\times$
provides a quasi-inverse
$(\underline M\text{-}\mathbf{Inv})^\times\to
\underline M^\times\text{-}\mathbf{Inv}$.

Especially, if $(X,\underline M)$ is a log scheme, we see
that the category of invertible $\underline M$-modules is
equivalent to that of invertible $\cO_{\!X}^\times$-modules,
hence also to that of invertible $\cO_{\!X}$-modules.

If $\phi:(Z,\underline N)\to(X,\underline M)$ is a morphism
of log schemes, and $\cN$ a $\underline M$-module, we let :
\set\begin{equation}\label{eq_pull-bak-mon-mnd}
\phi^*\!\cN:=\phi^{-1}\!\cN\otimes_{\phi^{-1}\underline M}\underline N.
\end{equation}
Clearly, $\phi^*\!\cN$ is an invertible $\underline N$-module,
whenever $\cN$ is an invertible $\underline M$-module.

\sset\subsubsection{}\label{subsec_log-proj}
Keep the notation of definition \ref{def_qcoh-log-gr}(i);
by faithfully flat descent, the restriction of $\cA$ to the
Zariski site of $X$ is again a quasi-coherent $\cO_{\!X}$-algebra,
which we denote again by $\cA$. We may then set $Y:=\Proj\,\cA$,
and let $\pi:Y\to X$ be the natural morphism. The composition of
$\pi^{-1}\beta_{\!\cA}$ and the morphism \eqref{eq_especiallly},
yields a map on the site $Y_\tau$ :
\set\begin{equation}\label{eq_yyield}
\pi^{-1}\underline P{}_n\to\cO_Y(n)
\qquad
\text{for every $n\in\N$}.
\end{equation}
Set $\cO_Y(\bullet):=\coprod_{n\in\Z}\cO_Y(n)$; it is
easily seen that the morphisms \eqref{eq_annother-onne}
induce a natural $\Z$-graded monoid structure on $\cO_Y(\bullet)$,
and the coproduct of the maps \eqref{eq_yyield} amounts
to a morphism of graded monoids :
$$
\omega_\bullet:\pi^{-1}\underline P\to\cO_Y(\bullet).
$$
We let $\underline Q$ be the push-out in the cocartesian diagram :
\set\begin{equation}\label{eq_grad-arrows}
{\diagram \omega_\bullet^{-1}(\cO_Y(\bullet)^\times) \ar[r] \ar[d] &
\pi^{-1}\underline P \ar[d] \\
\cO_Y(\bullet)^\times \ar[r] & \underline Q.
\enddiagram}
\end{equation}
Clearly $\underline Q$ is naturally a $\Z$-graded monoid,
in such a way that all the arrows in \eqref{eq_grad-arrows}
are morphisms of $\Z$-graded monoids. For every $n\in\Z$, let
$\underline Q{}_n$ be the degree $n$ subsheaf of $\underline Q$;
the map $\omega_\bullet$ and the natural inclusion
$\cO_Y(\bullet)^\times\to\cO_Y(\bullet)$ determine a unique
morphism $\underline Q\to\cO_Y(\bullet)$, whose restriction
in degree zero is a pre-log structure :
$$
\beta^\sim_{\!\cA}:\underline Q{}_0\to\cO_Y.
$$
Clearly $\alpha$ induces a unique morphism
$\alpha^\sim:\pi^{-1}\underline M\to\underline Q{}_0$,
such that the diagram of monoids :
$$
\xymatrix{ \pi^{-1}\underline M \ar[r]^-{\alpha^\sim}
\ar[d]_{\pi^{-1}\beta} &
\underline Q{}_0 \ar[d]^{\beta^\sim_{\!\cA}} \\
\pi^{-1}\underline\cO_{\!X} \ar[r]^-{\pi^\natural} & \cO_Y
}$$
commutes. Denote by $\underline P^\sim$ the log structure
associated to $\beta^\sim_{\!\cA}$; the
{\em homogeneous spectrum\/} of the quasi-coherent
$\N$-graded algebra $(\cA,\underline P,\alpha,\beta_{\!\cA})$
is defined as the $(X,\underline M)$-scheme :
$$
\Proj(\cA,\underline P):=(Y,\underline P^\sim).
$$
We also let :
$$
U_1(\cA,\underline P):=U_1(\cA)\times_Y\Proj(\cA,\underline P).
$$
Furthermore, for every $n\in\Z$ we have a natural morphism
of $\underline Q{}_0$-monoids :
\set\begin{equation}\label{eq_shift-mond}
\underline Q{}_0\otimes_{\cO^\times_Y}\cO_Y(n)^\times\to
\underline Q{}_n
\end{equation}
and it is easily seen that $\eqref{eq_shift-mond}_{|U_1(\cA)}$
is an isomorphism. We set :
$$
\underline P^\sim(n):=
(\underline P^\sim\otimes_{\underline Q{}_0}\underline Q{}_n)_{|U_1(\cA)}
\qquad\text{for every $n\in\Z$}.
$$
Hence \eqref{eq_shift-mond} induces a natural isomorphism :
\set\begin{equation}\label{eq_twist-of-logstr}
(\underline P^\sim\otimes_{\cO^\times_Y}\cO_Y(n)^\times)_{|U_1(\cA)}
\isom\underline P^\sim(n).
\end{equation}
Especially, $\underline P^\sim(n)$ is an invertible
$\underline P^\sim_{|U_1(\cA)}$-module, for every $n\in\Z$.
From \eqref{eq_twist-of-logstr}, we also deduce natural isomorphisms
of $\underline P^\sim_{|U_1(\cA)}$-modules :
\set\begin{equation}\label{eq_analog-never-dies}
\underline P^\sim(n)\otimes_{\underline P^\sim}\underline P^\sim(m)
\isom\underline P^\sim(n+m)
\qquad\text{for every $n,m\in\Z$}
\end{equation}
and of $\cO_{U_1(\cA)}$-modules :
\set\begin{equation}\label{eq_pass-to-O}
\underline P^\sim(n)\otimes_{\underline P^\sim}\cO_{U_1(\cA)}\isom
\cO_Y(n)_{|U_1(\cA)}
\qquad\text{for every $n\in\Z$}.
\end{equation}
Additionally, the morphism $\pi^{-1}\underline P{}_n\to\underline Q{}_n$
deduced from \eqref{eq_grad-arrows}, yields a natural map of
$\underline P^\sim_{|U_1(\cA)}$-modules :
\set\begin{equation}\label{eq_as_I_said}
\lambda_n:(\pi^*\underline P{}_n)_{|U_1(\cA)}\to\underline P^\sim(n)
\qquad\text{for every $n\in\N$.}
\end{equation}

\begin{example}\label{ex_invert-algebras}
Let $(Z,\gamma:\underline N\to\cO_{\!Z})$ be a log scheme,
$\cL$ an invertible $\underline N$-module, and set :
$$
\beta_{\cA(\cL)}:=
\Sym^\bullet_{\underline N}\cL\otimes_{\underline N}\gamma:
\Sym^\bullet_{\underline N}\cL\to\cA(\cL):=\Sym^\bullet_{\cO_{\!Z}}
(\cL\otimes_{\underline N}\cO_{\!Z})
$$
which is a morphisms of $\N$-graded monoids (notation of example
\ref{ex_symm-pow}). Clearly $\cA(\cL)$ is also a $\N$-graded
quasi-coherent $\cO_{\!Z}$-algebra. Denote also :
$$
\alpha_\cL:\underline N\to\Sym^\bullet_{\underline N}\cL
$$
the natural morphism that identifies $\underline N$ to
$\Sym^0_{\underline N}\cL$; then the datum
$$
(\cA(\cL),\Sym^\bullet_{\underline N}\cL,\alpha_\cL,\beta_{\cA(\cL)})
$$
is a quasi-coherent $\cO_{(Z,\underline N)}$-algebra, and
a direct inspection of the definitions shows that the induced
morphism of log schemes :
\set\begin{equation}\label{eq_proj-simple}
\pi_{(Z,\underline N)}:\P(\cL):=
\Proj(\cA(\cL),\Sym^\bullet_{\underline N}\cL)\to(Z,\underline N)
\end{equation}
is an isomorphism. Furthermore, we have natural isomorphisms
as in \eqref{eq_pulllback} :
$$
\pi_{(Z,\underline N)}^*(\cL^{\otimes n}\otimes_{\underline N}\cO_Z)
\isom\cO_{\P(\cL)}(n)
\qquad\text{for every $n\in\Z$}.
$$
Let $\underline P^\sim_\cL\to\cO_{\P(\cL)}$ be the log structure of
$\P(\cL)$; there follows a natural identification :
\set\begin{equation}\label{eq_ident-cL-and-twist}
\pi_{(Z,\underline N)}^*\cL^{\otimes n}\isom
\underline P^\sim_\cL\otimes_{\cO_{\P(\cL)}^\times}\cO_{\P(\cL)}(n)^\times
\isom\underline P^\sim_\cL(n)
\qquad\text{for every $n\in\Z$}
\end{equation}
where the last isomorphism is \eqref{eq_twist-of-logstr},
in view of the fact that $U_1(\cA(\cL))=\Proj\,\cA(\cL)$.
\end{example}

\begin{example}\label{ex_proj-space}
(i)\ \
Let $(Z,\underline N)$ be a log scheme, and $n\in\N$ any integer.
We define an $\N$-grading on $\N^{\oplus n}$, by setting
$$
\gr^k\N^{\oplus n}:=\{a_\bullet:=(a_1,\dots,a_n)~|~a_1+\cdots+a_n=k\}
\qquad
\text{for every $k\in\N$}.
$$
We then define the $\N$-graded monoid
$$
\Sym^\bullet_{\underline N}\underline N^{\oplus n}:=
\N^{\oplus n}_Z\times\underline N
$$
whose $\N$-grading is deduced in the obvious way from the
foregoing grading of $\N^{\oplus n}$.
The log structure $\gamma:\underline N\to\cO_{\!Z}$ extends
naturally to a map of $\N$-graded $Z$-monoids :
$$
\Sym_{\underline N}^\bullet\gamma^{\oplus n}:
\Sym^\bullet_{\underline N}\underline N^{\oplus n}\to
\Sym^\bullet_{\cO_{\!Z}}\cO_{\!Z}^{\oplus n}.
$$
Namely, if $e_1,\dots,e_n$ is the canonical basis of the
free $\cO_{\!Z}$-module $\cO_{\!Z}^{\oplus n}$, then
$\Sym^k_{\cO_{\!Z}}\cO_{\!Z}^{\oplus n}$ is a free
$\cO_{\!Z}$-module with basis
$$
\{e_{a_\bullet}:=
e^{a_1}_1\cdots e_n^{a_n}~|~a_\bullet\in\gr^k\N^{\oplus n}\}
$$
and $\Sym_{\underline N}^\bullet\gamma^{\oplus n}$ is given by
the rule : $(a_\bullet,x)\mapsto\gamma(x)\cdot e_{a_\bullet}$
for every $a_\bullet\in\N^{\oplus n}$, every $\tau$-open
subset of $Z$, and every section $x\in\underline N(U)$.
Clearly $\Sym_{\underline N}^\bullet\gamma^{\oplus n}$ defines
an $\N$-graded $\cO_{(Z,\underline N)}$-algebra
$$
\Sym^\bullet_{(\cO_{\!Z},\underline N)}
(\cO_{\!Z},\underline N)^{\oplus n}:=
(\Sym^\bullet_{\cO_{\!Z}}\cO_{\!Z}^{\oplus n},
\Sym^\bullet_{\underline N}\underline N^{\oplus n})
\qquad
\text{for every $n\in\N$}.
$$
We set :
$$
\P^n_{(Z,\underline N)}:=\Proj\,\Sym^\bullet_{(\cO_{\!Z},\underline N)}
(\cO_{\!Z},\underline N)^{\oplus n+1}
$$
and we call it the {\em projective $n$-dimensional space\/}
over $(Z,\underline N)$.

(ii)\ \
Denote by $\underline N^\sim$ the log structure of
$\P^n_{(Z,\underline N)}$, and by $\underline N^\sim(k)$ the
$\underline N^\sim$-modules defined as in \eqref{eq_twist-of-logstr},
for every $k\in\Z$. By simple inspection we get a commutative diagram
of monoids :
$$
{\diagram
\Gamma(Z,\Sym^k_{\underline N}\underline N^{\oplus n+1})
\ar[rrr]^-{\Gamma(Z,\Sym^k_{\underline N}\gamma^{\oplus n+1})} \ar[d]
& & & \Gamma(Z,\Sym^k_{\cO_{\!Z}}\cO_{\!Z}^{\oplus n+1}) \ar[d] \\
\Gamma(\P^n_{(Z,\underline N)},\underline N^\sim(k)) \ar[rrr] & & &
\Gamma(\P^n_{(Z,\underline N)},\cO_{\P^n_{(Z,\underline N)}}(k))
\enddiagram}
\qquad
\text{for every $k\in\N$}
$$
whose right vertical arrow is an isomorphism. Especially,
the natural basis of the free $\underline N(Z)$-module
$\Gamma(Z,\Sym^1_{\underline N}\underline N^{\oplus n+1})=
\underline N(Z)^{\oplus n+1}$ yields a distinguished system of
$n+1$ elements
$$
\eps_0,\dots,\eps_n\in
\Gamma(\P^n_{(Z,\underline N)},\underline N^\sim(1)).
$$
On the other hand, we have as well the distinguished system
of global sections
$$
T_0,\dots,T_n\in\Gamma(\P^n_{(Z,\underline N)},
\cO_{\P^n_{(Z,\underline N)}}(1))
$$
corresponding to the natural basis of
$\Gamma(Z,\Sym^1_{\cO_{\!Z}}\cO_{\!Z}^{\oplus n+1})=
\cO_{\!Z}(Z)^{\oplus n+1}$. For each $i=0,\dots,n$, the largest
open subset $U_i\subset\P^n_{(Z,\underline N)}$ such that
$T_i\in\Gamma(U_i,\cO_{\P^n_{(Z,\underline N)}}(1)^\times)$ is
the complement of the hyperplane where $T_i$ vanishes. Moreover,
notice that
$$
T_i^{-1}T_j\in\underline N^\sim(U_i)
\qquad
\text{for every $i,j=0,\dots,n$}.
$$
With this notation, the isomorphism \eqref{eq_twist-of-logstr}
yields the identification :
$$
\eps_j=T_i^{-1}T_j\otimes T_i
\qquad\text{on $U_i$}\qquad
\text{for every $i,j=0,\dots,n$}
$$
from which we also see that, for every $i=0,\dots,n$, the open
subset $U_i$ is the largest such that
$\eps_i\in\Gamma(U_i,\underline N^\sim(1)^\times)$. By the same
token, we obtain :
\set\begin{equation}\label{eq_trivial-proj-sp}
(\P^n_{(Z,\underline N)},\underline N^\sim)_\tr=U_0\cap\cdots\cap U_n
\simeq\G_{m,Z}^n.
\end{equation}
\end{example}

\sset\subsubsection{}\label{subsec_log-of-Gs}
Let $(Z,\underline N)$ be a log scheme,
$(\phi,\log\phi):(\cA,\underline P)\to(\cA',\underline P')$
a morphism of quasi-coherent $\N$-graded
$\cO_{(Z,\underline N)}$-algebras. We let (notation of
\eqref{subsec_glueee}) :
$$
G(\phi,\log\phi):=G(\phi)\times_{\Proj\,\cA'}\Proj(\cA',\underline P').
$$
Denote also by $\pi:Y:=\Proj\,\cA\to Z$ and $\pi':Y':=\Proj\,\cA'\to Z'$
the natural projections; there follows, on the one hand, a morphism
of $\N$-graded monoids :
\set\begin{equation}\label{eq_one.hand-mnd}
\pi^{\prime-1}(\log\phi):(\Proj\,\phi)^{-1}(\pi^{-1}\underline P)\to
(\pi^{\prime-1}\underline P')_{|G(\phi)}
\end{equation}
and on the other hand, a morphism of $\Z$-graded monoids :
\set\begin{equation}\label{eq_oth.hand-mnd}
(\Proj\,\phi)^{-1}\cO_Y(\bullet)^\times\to
\cO_{Y'}(\bullet)^\times_{|G(\phi)}
\end{equation}
deduced from \eqref{eq_new-nu}. Define the $\Z$-graded monoid
$\underline Q$ on $Y_\tau$ as in \eqref{eq_grad-arrows}, and
the analogous $\Z$-graded monoid $\underline Q'$ on $Y'_\tau$.
Then \eqref{eq_one.hand-mnd} and \eqref{eq_oth.hand-mnd} determine
a unique morphism of $\Z$-graded monoids :
$$
\theta:(\Proj\,\phi)^{-1}\underline Q\to\underline Q'_{|G(\phi)}
$$
and by construction, the restriction of $\theta$ in degree zero
is a morphism of pre-log structures :
$$
\theta_0:(\Proj\,\phi)^*((\underline Q)_0,\beta^\sim_{\!\cA})\to
((\underline Q')_0,\beta^\sim_{\!\cA'})
$$
whence a morphism of $(Z,\underline N)$-schemes :
$$
\Proj(\phi,\log\phi):G(\phi,\log\phi)\to\Proj(\cA,\underline P).
$$
Moreover, on the one hand, \eqref{eq_new-nu} induces an isomorphism
of $\cO^\times_{Y'}$-modules :
$$
\nu_{\!\cA(n)}^\times:
(\cO^\times_{Y'}\otimes_{(\Proj\,\phi)^{-1}\cO^\times_Y}
(\Proj\,\phi)^{-1}\cO_Y(n)^\times)_{|G_1(\phi)}\isom
\cO_{Y'}(n)^\times_{|G_1(\phi)}
\qquad\text{for every $n\in\Z$}
$$
(notation of \eqref{eq_caveatt}). On the other hand, for every
$n\in\Z$, the morphism of
$(\Proj\,\phi)^{-1}(\underline Q)_0$-modules $\theta_n$
determines a morphism of $\underline P^{\prime\sim}$-modules :
\set\begin{equation}\label{eq_exhausting-theta}
\theta^\sim_n:\Proj(\phi,\log\,\phi)^*\underline P^\sim(n)_{|G_1(\phi)}
\to\underline P^{\prime\sim}(n)_{|G_1(\phi)}
\end{equation}
and by inspecting the construction, it is easily seen that
the isomorphism \eqref{eq_twist-of-logstr} (and the corresponding
one for $\underline P^{\prime\sim}(n)$) identifies $\theta^\sim_n$
with $\nu_{\!\cA(n)}^\times\otimes_{\cO^\times_{Y'}}
\underline P^{\prime\sim}$; especially, $\theta^\sim_n$ is
an isomorphism.

\sset\subsubsection{}
Let $\psi:(Z',\underline N')\to(Z,\underline N)$ be a morphism
of log schemes, and $(\cA,\underline P,\alpha,\beta_\cA)$
a $\N$-graded quasi-coherent $\cO_{(Z,\underline N)}$-algebra.
We may view $\underline P$ as a $\underline N$-module, via the
morphism $\alpha$, hence we may form the $\underline N'$-module
$\psi^*\underline P$, as in \eqref{eq_pull-bak-mon-mnd}.
Moreover, by remark \eqref{rem_tens-is-pushout}(i),
$\psi^*\underline P$ is a $\N$-graded sheaf of monoids
on $Z'_\tau$, such that
\set\begin{equation}\label{eq_pulbak-gr-algs}
\psi^*(\cA,\underline P):=
(\psi^*\cA,\psi^*\underline P,\psi^*\alpha,\psi^*\beta_\cA)
\end{equation}
is a $\N$-graded quasi-coherent $\cO_{(Z',\underline N')}$-algebra,
and in view of the isomorphism \eqref{eq_tens-is-pushout}, we obtain
a natural isomorphism of $(Z',\underline N')$-schemes:
\set\begin{equation}\label{eq_proj-of-puback}
\Proj\,\psi^*(\cA,\underline P)\isom
(Z',\underline N')\times_{(Z,\underline N)}\Proj(\cA,\underline P).
\end{equation}
Furthermore, denote by $\pi_{(\cA,\underline P)}:
U_1(\psi^*(\cA,\underline P))\to U_1(\cA,\underline P)$
the morphism deduced from \eqref{eq_proj-of-puback},
and by $\pi_Y:Y':=\Proj\,\psi^*\cA\to Y:=\Proj\,\cA$ the
underlying morphism of schemes. From \eqref{eq_for-arbi-n}
we obtain natural isomorphisms :
$$
(\cO_{Y'}^\times\otimes_{\pi_Y^{-1}\cO^\times_Y}
\pi_Y^{-1}\cO_Y(n)^\times)_{|U_1(\psi^*\cA)}\isom
\cO_{Y'}(n)^\times_{|U_1(\psi^*\cA)}
\qquad\text{for every $n\in\Z$}
$$
and the latter induce natural identifications :
\set\begin{equation}\label{eq_last-small-one}
\pi^*_{(\cA,\underline P)}\underline P^\sim(n)\isom
(\psi^*\underline P)^\sim(n)
\qquad\text{for every $n\in\Z$}.
\end{equation}

\sset\subsubsection{}
Keep the notation of \eqref{subsec_log-proj}, and let
$\log\cC_{(X,\underline M)}$ be the category whose objects are
the pairs $((Z,\underline N),\cL)$, where $(Z,\underline N)$
is a $(X,\underline M)$-scheme, and $\cL$ is an invertible
$\underline N$-module. The morphisms
$((Z,\underline N),\cL)\to((Z',\underline N'),\cL')$
are the pairs $(\phi,h)$, where
$\phi:(Z,\underline N)\to(Z',\underline N')$ is a morphism
of $(X,\underline M)$-schemes, and $h:\phi^*\cL'\isom\cL$
is an isomorphism of $\underline N$-modules (with composition
of morphisms defined in the obvious way). There is an obvious
forgetful functor :
$$
\sp:\log\cC_{(X,\underline M)}\to\cC_X
\quad : \quad
((Z,\underline N),\cL)\mapsto(Z,\cL\otimes_{\underline N}\cO_{\!Z})
$$
and the functor $F_{\!\cA}$ can be lifted to a functor :
$$
F_{(\cA,\underline P)}:\log\cC_{(X,\underline M)}\to\Set
$$
which assigns to any object
$((Z,\gamma:\underline N\to\cO_{\!Z}),\cL)$ the set consisting
of all morphisms of $\N$-graded quasi-coherent
$\cO_{(Z,\underline N)}$-algebras :
$$
\xymatrix{
\psi^*\underline P \ar[rr]^-{\log g} \ar[d]_{\psi^*\beta_\cA}
& & \Sym^\bullet_{\underline N}\cL \ar[d]^{\beta_{\cA(\cL)}} \\
\psi^*\cA \ar[rr]^-g & & \cA(\cL)
}$$
where $\psi:(Z,\underline N)\to(X,\underline M)$ is the
structural morphism, and $g$ is an epimorphism on the underlying
$\cO_{\!Z}$-modules.

\begin{proposition}\label{prop_represent-log}
$(U_1(\cA,\underline P),\underline P^\sim(1))\in
\Ob(\log\cC_{(X,\underline M)})$ represents the functor
$F_{(\cA,\underline P)}$.
\end{proposition}
\begin{proof} Given
$((Z,\underline N)\xrightarrow{\psi}(X,\underline M),\cL)\in
\Ob(\log\cC_{(X,\underline M)})$ and
$(g,\log g)\in F_{(\cA,\underline P)}((Z,\underline N),\cL)$,
define $\P(\cL)$ as in \eqref{eq_proj-simple}; there follows
a morphism of $(Z,\underline N)$-schemes :
$$
\Proj(g,\log g):\P(\cL)\to\Proj\,\psi^*(\cA,\underline P).
$$
In view of \eqref{eq_proj-simple} and \eqref{eq_proj-of-puback},
this is the same as a morphism of $(X,\underline M)$-schemes :
$$
\P(g,\log g):(Z,\underline N)\to\Proj(\cA,\underline P)
$$
and arguing as in the proof of lemma \ref{lem_repres-proj},
we see that the image of $\P(g,\log g)$ lands in
$U_1(\cA,\underline P)$. Next, combining
\eqref{eq_ident-cL-and-twist}, \eqref{eq_exhausting-theta}
and \eqref{eq_last-small-one}, we deduce a natural isomorphism :
$$
\begin{aligned}
\pi_{(Z,\underline N)}^*\circ\P(g,\log g)^*\underline P^\sim(1)
& \isom\Proj(g,\log g)^*\circ
\pi_{(\cA,\underline P)}^*\underline P^\sim(1) \\
& \isom\Proj(g,\log g)^*(\psi^*\underline P)^\sim(1) \\
& \isom\underline P_\cL^\sim(1) \\
& \isom\pi_{(Z,\underline N)}^*\cL
\end{aligned}
$$
whence an isomorphism
$h(g,\log g):\P(g,\log g)^*\underline P^\sim(1)\isom\cL$,
and the datum :
$$
(\P(g,\log g),h(g,\log g)):((Z,\underline N),\cL)\to
(U_1(\cA,\underline P),\underline P^\sim(1))
$$
is a well defined morphism of $\log\cC_{(X,\underline M)}$.

Conversely, let $\phi:=(\beta,\log\beta):
(Z,\underline N)\to U_1(\cA,\underline P)$ be a morphism
of $(X,\underline M)$-schemes, and
$h:\phi^*\underline P^\sim(1)\isom\cL$ an isomorphism of
$\underline N$-modules. In view of \eqref{eq_analog-never-dies},
we deduce an isomorphism :
$$
h^{\otimes n}:\phi^*\underline P^\sim(n)\isom\cL^{\otimes n}
\qquad\text{for every $n\in\Z$}.
$$
Combining with \eqref{eq_as_I_said}, we may define the map of
$\underline N$-modules :
$$
\log\hat g(\phi,h):=
\bigoplus_{n\in\N}h^{\otimes n}\circ\beta^*(\lambda_n):
\psi^*\underline P\to\Sym_{\underline N}^\bullet\cL.
$$
On the other hand, in view of \eqref{eq_pass-to-O}, we have an
isomorphism of $\cO_{\!Z}$-modules :
$$
h\otimes_{\underline N}\cO_{\!Z}:\beta^*\cO_Y(1)_{|U_1(\cA)}\isom
\phi^*\underline P^\sim(1)\otimes_{\underline N}\cO_{\!Z}\isom
\cL\otimes_{\underline N}\cO_{\!Z}
$$
(where, as usual, $Y:=\Proj\,\cA$). We let (notation of
\eqref{eq_promise}) :
$$
\hat g(\phi,h):=g(\beta,h\otimes_{\underline N}\cO_{\!Z})
$$
and notice that the pair $(\hat g(\phi,h),\log\hat g(\phi,h))$
is an element of $F_{(\cA,\underline P)}((Z,\underline N),\cL)$.
Summing up, we have exhibited two natural transformations :
$$
\begin{aligned}
\theta:F_{(\cA,\underline P)}\Rightarrow
\Hom_{\log\cC_{(X,\underline M)}}
(-,(U_1(\cA,\underline P),\underline P^\sim(1)))
& & (g,\log g) \mapsto\, & \P(g,\log g) \\
\sigma:\Hom_{\log\cC_{(X,\underline M)}}
(-,(U_1(\cA,\underline P),\underline P^\sim(1)))
\Rightarrow F_{(\cA,\underline P)}
& & (\phi,h) \mapsto\, & (\hat g(\phi,h),\log\hat g(\phi,h))
\end{aligned}
$$
and it remains to show that these transformations are isomorphisms
of functors. However, the latter fit into an essentially commutative
diagram of natural transformations :
$$
\xymatrix{
F_{(\cA,\underline P)} \ar@{=>}[r]^-\theta \ar@{=>}[d] &
\Hom_{\log\cC_{(X,\underline M)}}
(-,(U_1(\cA,\underline P),\underline P^\sim(1)))
 \ar@{=>}[r]^-\sigma \ar@{=>}[d] & F_{(\cA,\underline P)} \ar@{=>}[d] \\
F_{\!\cA}\circ\sp \ar@{=>}[r] &
\Hom_{\cC_X}(\sp(-),(U_1(\cA),\cO_Y(1)_{|U_1(\cA)})) \ar@{=>}[r]
& F_{\!\cA}\circ\sp
}$$
whose bottom line is given by the natural transformations
\eqref{eq_rules}. Moreover, given an isomorphism $h$ of invertible
$\underline N$-modules, the discussion in \eqref{subsec_inv-mod-mon}
leads to the identity :
\set\begin{equation}\label{eq_leadstothe}
h=(h\otimes_{\underline N}\cO_{\!Z})^\times
\otimes_{\cO_{\!Z^\times}}\underline N.
\end{equation}
Likewise, we have natural identifications :
$$
\Sym^n_{\underline N}\cL=
(\Sym^n_{\cO_{\!Z}}\cL\otimes_{\underline N}\cO_{\!Z})^\times
\otimes_{\cO_{\!Z}^\times}\underline N
\qquad\text{for every $n\in\Z$}
$$
which show that $\beta_{\!\cA}$ and $g$ determine uniquely $\log g$.
This -- and an inspection of the proof of lemma \ref{lem_repres-proj} --
already implies that $\sigma\circ\theta$ is the identity automorphism
of the functor $F_{(\cA,\underline P)}$. Finally, let
$\phi=(\beta,\log\beta)$ and $h$ as in the foregoing,
so that $(\phi,h)$ is a morphism in $\log\cC_{(X,\underline M)}$; to
conclude we have to show that $(\phi',h'):=\theta\circ\sigma(\phi,h)$
equals $(\phi,h)$. Say that $\phi':=(\beta',\log\beta')$;
by the above (and by the proof of lemma \ref{lem_repres-proj})
we already know that $\beta=\beta'$, and \eqref{eq_leadstothe}
implies that $h=h'$. Hence, it remains only to show that
$\log\beta=\log\beta'$, which can be checked directly on
the stalks over the $\tau$-points of $Z$ : we leave
the details to the reader.
\end{proof}

\begin{example}\label{ex_proj-space-continued}
(i)\ \ 
Let $\psi:(Z',\underline N')\to(Z,\underline N)$ be a morphism of
log schemes, $n\in\N$ any integer, and $\P^n_{(Z,\underline N)}$
the $n$-dimensional projective space over $(Z,\underline N)$,
defined as in example \ref{ex_proj-space}. A simple inspection
of the definitions yields a natural isomorphism of
$\cO_{(Z',\underline N')}$-algebras
$$
\Sym^\bullet_{(\cO_{\!Z'},\underline N')}
(\cO_{\!Z'},\underline N')^{\oplus n}\isom
\psi^*\Sym^\bullet_{(\cO_{\!Z},\underline N)}
(\cO_{\!Z},\underline N)^{\oplus n}
\qquad
\text{for every $n\in\N$}
$$
whence a natural isomorphism of $(Z',\underline N')$-schemes :
$$
\P^n_{(Z',\underline N')}\isom
(Z',\underline N')\times_{(Z,\underline N)}\P^n_{(Z,\underline N)}
\qquad
\text{for every $n\in\N$}.
$$

(ii)\ \
Let $\cL$ be any invertible $\underline N$-module; notice that
a morphism of $\N$-graded monoids
$$
\Sym^\bullet_{\underline N}\underline N^{\oplus n}\to
\Sym^\bullet_{\underline N}\cL
$$
which is the identity map in degree zero, is the same as the
datum of a sequence
$$
(\beta_i:\underline N\to\cL~|~i=1,,\dots,n)
$$
of morphisms of $\underline N$-modules, and the latter is the
same as a sequence $(b_1,\dots,b_n)$ of global sections of
$\cL$. Since $\Sym^1_{\cO_{\!Z}}\cO_{\!Z}^{\oplus n+1}$
generates $\Sym^\bullet_{\cO_{\!Z}}\cO_{\!Z}^{n+1}$, proposition
\ref{prop_represent-log} and (i) imply that $\P^n_{(Z,\underline N)}$
represents the functor $\log\cC_{(Z,\underline N)}\to\Set$ that
assigns to any pair $((X,\underline M),\cL)$ the set of all
sequences $(b_0,\dots,b_n)$ of global sections of $\cL$. The
bijection
\set\begin{equation}\label{eq_many-parentheses}
\Hom_{\log\cC_{(Z,\underline N)}}
(((X,\underline M),\cL),(\P^n_{(Z,\underline N)},\underline N^\sim(1)))
\isom\Gamma(X,\cL)^{n+1}
\end{equation}
can be explicited as follows. Let
$(\phi,h):((X,\underline M),\cL)\to
(\P^n_{(Z,\underline N)},\underline N^\sim(1))$ be a given morphism
in $\log\cC_{(Z,\underline N)}$; then
$h:\phi^*\underline N^\sim(1)\to\cL$ is an isomorphism of
$\underline M$-modules, which induces a map on global sections
$$
\Gamma(h):\Gamma(\P^n_{(Z,\underline N)},\underline N^\sim(1)))\to
\Gamma(X,\cL).
$$
Now, $\underline N^\sim(1)$ admits a distiguished system of global
sections $\eps_0,\dots,\eps_n$ (example \ref{ex_proj-space}(ii)),
and the bijection \eqref{eq_many-parentheses} assigns to $(\phi,h)$
the sequence $(\Gamma(h)(\eps_0),\dots,\Gamma(h)(\eps_n))$.

(iii)\ \
Given a sequence $b_\bullet:=(b_0,\dots,b_n)$ as in (ii),
denote by
$$
f_{b_\bullet}:(X,\underline M)\to\P^n_{(Z,\underline N)}
$$
the corresponding morphism. A direct inspection of the definitions
shows that the formation of $f_{b_\bullet}$ is compatible with
arbitrary base changes $h:(Z',\underline N')\to(Z,\underline N)$.
Namely, set $(X',\underline M'):=
(Z',\underline N')\times_{(Z,\underline N)}(X,\underline M)$,
let $g:(X',\underline M')\to(X,\underline M)$ be the induced
morphism, $\cL':=g^*\cL$, and suppose that $\cL'$ is also
invertible (which always holds, if $\underline M'$ is an integral
log structure); the sequence $b_\bullet$ pulls back to a corresponding
sequence $b'_\bullet:=(b'_0,\dots,b_n)$ of global sections of $\cL'$,
and there follows a cartesian diagram of log schemes :
$$
\xymatrix{
(X',\underline M') \ar[r]^-{f_{b'_\bullet}} \ar[d]_g &
\P^n_{(Z',\underline N')} \ar[d]^{\P^n_h} \\
(X,\underline M) \ar[r]^-{f_{b_\bullet}} & \P^n_{(Z,\underline N)}.
}$$
The details shall be left to the reader.
\end{example}

\begin{definition}\label{def_blow-up}
Let $(X,\underline M)$ be a log scheme, and
$\cI\subset\underline M$ an ideal of $\underline M$
(see \eqref{subsec_sheaves-of-mons}).
\begin{enumerate}
\item
Let $f:(Y,\underline N)\to(X,\underline M)$ be a morphism
of log schemes; then $f^{-1}\!\cI$ is an ideal in the sheaf
of monoids $f^{-1}\underline M$; we let :
$$
\cI\underline N:=\log f(f^{-1}\!\cI)\cdot\underline N
$$
which is the smallest ideal of $\underline N$ containing
the image of $f^{-1}\!\cI$.
\item
A {\em logarithmic blow up\/} of the ideal $\cI$ is a
morphism of log schemes
$$
\phi:(X',\underline M')\to(X,\underline M)
$$
which enjoys the following universal property. The ideal
$\cI\underline M'$ is invertible, and every morphism
of log schemes $(Y,\underline N)\to(X,\underline M)$
such that $\cI\underline N\subset\underline N$ is invertible,
factors uniquely through $\phi$.
\end{enumerate}
\end{definition}

\begin{remark}\label{rem_quite}
(i)\ \ Keep the notation of definition \ref{def_blow-up}.
By the usual general nonsense, it is clear that the blow
up $(X',\underline M')$ is determined up to unique isomorphism
of $(X,\underline M)$-schemes.

(ii)\ \
Moreover, let $f:Y\to X$ be a morphism of schemes. Then we claim
that the natural projection :
$$
(Y',\underline M'_Y):=Y\times_X(X',\underline M')\to
(Y,\underline M_Y):=Y\times_X(X,\underline M)
$$
is a logarithmic blow up of $\cI\underline M_Y$, provided
$\cI\underline M'_Y$ is an invertible ideal; especially, this holds
whenever $\underline M'$ is an integral log structure (lemma
\ref{lem_simple-charts-top}(i)). The proof is left as an exercise
for the reader.
\end{remark}

\sset\subsubsection{}\label{subsec_blow-up-alg}
Let $(X,\underline M)$ be a log scheme, $\cI\subset\underline M$
an ideal; we shall show the existence of the logarithmic blow up
of $\cI$, under fairly general conditions. To this aim, we introduce
the graded {\em blow up $\cO_{\!X}$-algebra} :
$$
\sB(X,\underline M,\cJ):=
\bigoplus_{n\in\N}\cI^n\otimes_{\underline M}\cO_{\!X}
$$
where $\cI^n\otimes_{\underline M}\cO_{\!X}$ is the sheaf associated
to the presheaf $U\mapsto\cI^n(U)\otimes_{\underline
M(U)}\cO_{\!X}(U)$ on $X_\tau$. Here $\cI^0:=\underline M$, and for
every $n>0$ we let $\cI^n$ be the sheaf associated to the presheaf
$U\mapsto\cI(U)^n$ on $X_\tau$. The graded multiplication law of the
blow up $\cO_{\!X}$-algebra is induced by the multiplication
$\cI^n\times\cI^m\to\cI^{n+m}$, for every $n,m\in\N$.

The natural map $\cI^n\to\cI^n\otimes_{\underline M}\cO_{\!X}$
induces a morphism of sheaves of monoids :
$$
\coprod_{n\in\N}\cI^n\to\sB(X,\underline P,\cJ).
$$
The latter defines a $\N$-graded $\cO_{(X,\underline M)}$-algebra,
which we denote $\cB(X,\underline M,\cI)$.

\sset\subsubsection{}\label{subsec_first-invert}
Suppose first that $\cI$ is invertible; then it is easily seen
that, locally on $X_\tau$, $\cI$ is generated by a regular local
section (see example \ref{ex_regular}(i)), hence the same holds for
the power $\cI^n$, for every $n\in\N$. Therefore $\cI^n$ is a locally
free $\underline M$-module of rank one, and we have a natural
isomorphism :
$$
\cB(X,\underline M,\cI)\isom(\cA(\cI),\Sym^\bullet_{\underline M}\cI)
$$
(notation of example \eqref{ex_invert-algebras}). It follows that in
this case, the natural projection :
$$
\pi_{(X,\underline M,\cI)}:\Proj\,\cB(X,\underline M,\cI)\to(X,\underline M)
$$
is an isomorphism of log schemes.

\sset\subsubsection{}\label{subsec_formation-blow}
The formation of $\cB(X,\underline M,\cI)$ is obviously functorial
with respect to morphisms of log schemes; more precisely,
such a morphism $f:(Y,\underline N)\to(X,\underline M)$
induces a morphism of $\N$-graded $\cO_{(Y,\underline N)}$-algebras:
$$
\cB(f,\cI):f^*\cB(X,\underline M,\cI)\to\cB(Y,\underline N,\cI\underline N)
$$
(notation of \eqref{eq_pulbak-gr-algs}) which is an epimorphism
on the underlying $\cO_Y$-modules.
Moreover, if $g:(Z,\underline Q)\to(Y,\underline N)$ is another
morphism, we have the identity :
\set\begin{equation}\label{eq_little-cocycle}
\cB(f\circ g,\cI)=\cB(g,\cI\underline N)\circ g^*\cB(f,\cI).
\end{equation}
Furthermore, the construction of the blow up algebra is
local for the topology of $X_\tau$ : if $U\to X$ is any
object of $X_\tau$, we have a natural identification
$$
\cB(U,\underline M_{|U},\cI_{|U})\isom\cB(X,\underline M,\cI)_{|U}.
$$
In the presence of charts for the log structure $\underline M$,
we can give a handier description for the blow up algebra;
namely, we have the following :

\begin{lemma}\label{lem_blow-alg}
Let $X$ be a scheme, $\alpha:\underline P\to\cO_{\!X}$ a
pre-log structure, $\beta:\underline P\to\underline P^{\log}$
the natural morphism of pre-log structures. Let also
$\cI\subset\underline P$ be an ideal, and set
$\cI\underline P^{\log}:=\beta(\cI)\cdot\underline P^{\log}$
(which is the ideal of $\underline P^{\log}$ generated
by the image of $\cI$). Then :
\begin{enumerate}
\item
There is a natural isomorphism of graded $\cO_{\!X}$-algebras :
$$
\bigoplus_{n\in\N}\cI^n\otimes_{\underline P}\cO_{\!X}\isom
\sB(X,\underline P^{\log},\cI\underline P^{\log}).
$$
\item
Especially, suppose $(X,\underline M)$ is a log scheme that admits a
chart $\beta:P_X\to\underline M$, let $I\subset P$ be an ideal, and
set $I\underline M:=\beta(I_X)\cdot\underline M$. Then
$\cB(X,\underline M,I\underline M)$ is a $\N$-graded quasi-coherent
$\cO_{(X,\underline M)}$-algebra.
\item
In the situation of\/ {\em (ii)}, set $(S,P^{\log}_S):=\Spec(\Z,P)$
(see \eqref{subsec_constant-log}), and denote by
$f:(X,\underline M)\to(S,P^{\log}_S)$ the natural morphism.
Then the map :
$$
\cB(f,IP^{\log}_S):f^*\cB(S,P^{\log}_S,IP^{\log}_S)
\to\cB(X,\underline M,I\underline M)
$$
is an isomorphism of $\N$-graded $\cO_{(X,\underline M)}$-algebras.
\end{enumerate}
\end{lemma}
\begin{proof}(i): Since the functor \eqref{eq_another-one} commutes
with all colimits, we have a natural isomorphism of sheaves of
rings on $X_\tau$ :
$$
\Z[\underline P^{\log}]\isom
\Z[\underline P]\otimes_{\Z[\alpha^{-1}\cO^\times_{\!X}]}\Z[\cO^\times_{\!X}].
$$
We are therefore reduced to showing that the natural morphism
$$
\Z[\cI^n]\otimes_{\Z[\alpha^{-1}\cO^\times_{\!X}]}\Z[\cO^\times_{\!X}]
\to\cI^n\Z[\underline P^{\log}]
$$
is an isomorphism for every $n\in\N$. The latter assertion
can be checked on the stalks over the $\tau$-points of $X$;
to this aim, we invoke the more general :

\begin{claim}\label{cl_invok-mr-gen}
Let $G$ be an abelian group, $\phi:H\to P'$ and $\psi:H\to G$
two morphisms of monoids, $I\subset P'$ an ideal. Then the
natural map
\set\begin{equation}\label{eq_blow-ideal}
\Z[I]\otimes_{\Z[H]}\Z[G]\to I\Z[P'\otimes_HG]
\end{equation}
is an isomorphism.
\end{claim}
\begin{pfclaim} Recall that $\Z[I]=I\Z[P']$, and the map
\eqref{eq_blow-ideal} is induced by the natural identification:
$\Z[P']\otimes_{\Z[H]}\Z[G]\isom\Z[P'\otimes_HG]$ (see
\eqref{eq_push-out-tensor}). Especially, \eqref{eq_blow-ideal}
is clearly surjective, and it remains to show that it is
also injective. To this aim, notice first that $\psi$ factors
through the unit of adjunction  $\eta:H\to H^\gp$; the morphism
$\Z[\eta]:\Z[H]\to\Z[H^\gp]=H^{-1}\Z[H]$ is a localization
map (see \eqref{eq_loc-monds}), especially it is flat, hence
\eqref{eq_blow-ideal} is injective when $G=H^\gp$ and $\psi=\eta$.
It follows easily that we may replace $H$ by $H^\gp$, $P'$
by $P'\otimes_HH^\gp$, $I$ by $I\cdot(P'\otimes_HH^\gp)$,
and therefore assume that $H$ is a group. Let $L:=\psi(H)$;
arguing as in the foregoing, we may consider separately the
case where $\psi$ is the surjection $H\to L$ and the case
where $\psi$ is the injection $i:L\to G$. However, one sees
easily that $\Z[i]:\Z[L]\to\Z[G]$ is a flat morphism, hence
it suffices to consider the case where $\psi$ is a surjective
group homomorphism. Set $K:=\Ker\,\psi$; we have a natural
identification : $P'\otimes_HG\simeq P'\otimes_K\{1\}$, hence
we may further reduce to the case where $G=\{1\}$. Then
the contention is that the augmentation map $\Z[K]\to\Z$
induces an isomorphism $\omega:\Z[I]\otimes_{\Z[K]}\Z\isom I\Z[P'/K]$.
From lemma \ref{lem_special-p-out}(ii), we derive easily that
$I\Z[P'/K]=\Z[I/K]$, where $I/K$ is the set-theoretic quotient
of $I$ for the $K$-action deduced from $\phi$. However,
any set-theoretic section $I/K\to I$ of the natural projection
$I\to I/K$ yields a well-defined surjection
$\Z[I/K]\to\Z[I]\otimes_{\Z[K]}\Z$ whose composition with
$\omega$ is the identity map. The claim follows.
\end{pfclaim}

(ii): Let $U$ be an affine object of $X_\tau$, say
$U=\Spec\,A$; from (i), we see that
$\cB(X,\underline M,I\underline M)_{|U}$ is the quasi-coherent
$\cO_{\!U}$-algebra associated to the $A$-algebra
$\oplus_{n\in\N}\Z[I^n]\otimes_{\Z[P]}A$.

(iii): In view of (i) we know already that $\cB(f,IP^{\log}_S)$
induces an isomorphism on the underlying $\cO_{\!X}$-algebras;
hence, by lemma \ref{lem_time}(i), it remains to show that
$\cB(f,IP^{\log}_S)$ induces an isomorphism :
$$
f^*(I^nP^{\log}_S)\isom I^nf^*(P^{\log}_S)
\qquad\text{for every $n\in\N$}.
$$
Let
$\gamma:f^{-1}P^{\log}_S\to f^{-1}\cO_S\to\cO_{\!X}$ be the
natural map; after replacing $I$ by $I^n$, we come down to
showing that the natural map :
$$
f^{-1}(IP^{\log}_S)\otimes_{\gamma^{-1}\cO^\times_{\!X}}
\cO^\times_{\!X}\to I\cdot
(f^{-1}P^{\log}_S\otimes_{\gamma^{-1}\cO^\times_{\!X}}\cO^\times_{\!X})
$$
is an isomorphism. This assertion can be checked on the stalks
over the $\tau$-points of $X$; if $x$ is such a point, let
$G:=\cO^\times_{\!X,x}$, $H:=\gamma^{-1}_x G$ and
$P':=P^{\log}_{S,f(x)}$. The map under consideration is
the natural morphism of $P'$-modules :
$$
\omega:(IP')\otimes_HG\to I(P'\otimes_HG)
$$
and it suffices to show that $\Z[\omega]$ is an isomorphism;
however, the latter is none else than \eqref{eq_blow-ideal},
so we may appeal to claim \ref{cl_invok-mr-gen} to conclude.
\end{proof}

\sset\subsubsection{}\label{subsec_coh-ideals-blow}
We wish to generalize lemma \ref{lem_blow-alg}(ii) to log
structures that do not necessarily admit global charts. Namely,
suppose now that $\underline M$ is a quasi-coherent log structure
on $X$, and $\cI\subset\underline M$ is a coherent ideal (see
definition \ref{def_coh-idea-log}(v)). For every $\tau$-point
$\xi$ of $X$, we may then find a neighborhood $U$ of $\xi$ in
$X_\tau$, a chart $\beta:P_U\to\underline M_{|U}$, and local
sections $s_1,\dots,s_n\in\cI(U)$ which form a system of generators
for $\cI_{|U}$. We may then write $s_{i,\xi}=u_i\cdot\beta(x_i)$ for
certain $x_1,\dots,x_n\in P$ and
$u_1,\dots,u_n\in\cO^\times_{\!X,\xi}$. Up to shrinking $U$, we may
assume that $u_1,\dots,u_n\in\cO^\times_{\!X}(U)$, and it follows
that $\cI_{|U}=I\cdot\underline M{}_{|U}$, where $I\subset P$ is the
ideal generated by $x_1,\dots,x_n$. In other words, locally on
$X_\tau$, the datum $(X,\underline M,\cI)$ is of the type considered
in lemma \ref{lem_blow-alg}(ii); therefore the blow up
$\cO_{\!X}$-algebra $\cB(X,\underline M,\cI)$ is quasi-coherent. We
may then consider the natural projection :
\set\begin{equation}\label{eq_blow-up-morph}
\pi_{(X,\underline M,\cI)}:\Bl_\cI(X,\underline M):=
\Proj\,\cB(X,\underline M,\cI)\to(X,\underline M).
\end{equation}
Next, let $f:(Y,\underline N)\to(X,\underline M)$
be a morphism of log schemes; it is easily seen that
$\cI\underline N$ is a coherent ideal of $\underline N$,
hence $\cB(Y,\underline N,\cI\underline N)$ is quasi-coherent
as well, and since the map $\cB(f,\cI)$ of \eqref{subsec_formation-blow}
is an epimorphism on the underlying $\cO_Y$-modules, we have :
$$
G(\cB(f,\cI))=\Bl_{\cI\underline N}(Y,\underline N)
$$
(notation of \eqref{subsec_log-of-Gs}) whence a closed immersion of
$(Y,\underline N)$-schemes :
$$
\Proj\,\cB(f,\cI):\Bl_{\cI\underline N}(Y,\underline N)\to
(Y,\underline N)\times_{(X,\underline M)}\Bl_\cI(X,\underline M)
$$
(\cite[Ch.II, Prop.3.6.2(i)]{EGAII} and \eqref{eq_proj-of-puback}),
which is the same as a morphism of $(X,\underline M)$-schemes :
$$
\phi(f,\cI):\Bl_{\cI\underline N}(Y,\underline N)\to
\Bl_\cI(X,\underline M).
$$
Moreover, if $g:(Z,\underline Q)\to(Y,\underline N)$ is another
morphism, \eqref{eq_little-cocycle} induces the identity :
\set\begin{equation}\label{eq_little-cocylce-blow}
\phi(f,\cI)\circ\phi(g,\cI\underline N)=\phi(f\circ g,\cI).
\end{equation}

\begin{example}\label{ex_explicit-blowup}
In the situation of lemma \ref{lem_blow-alg}(iii), notice that
$\cO_{\!S}$ is a flat $\Z[P_S]$-algebra, and therefore :
$$
\sB(S,P^{\log}_S,IP^{\log}_S)=\bigoplus_{n\in\N}I^n\cO_{\!S}.
$$
Moreover, \eqref{eq_proj-of-puback} specializes to a natural
isomorphism of $(X,\underline M)$-schemes :
\set\begin{equation}\label{eq_special-isoblow}
\Bl_{I\underline M}(X,\underline M)\isom
X\times_S\Bl_{IP^{\log}_S}(S, P^{\log}_S).
\end{equation}
We wish to give a more explicit description of the log
structure of $\Bl_{IP^{\log}_S}(S, P^{\log}_S)$.
To begin with, recall that $S':=\Proj\,\sB(S,P^{\log}_S,IP^{\log}_S)$
admits a distinguished covering by (Zariski) affine open subsets :
namely, for every $a\in I$, consider the localization
$$
P_a:=T^{-1}_aP
\qquad\text{where}\qquad
T_a:=\{a^n~|~n\in\N\}
$$
and let $Q_a\subset P_a$ be the submonoid generated by the image
of $P$ and the subset $\{a^{-1}b~|~b\in I\}$; then
$$
S'=\bigcup_{a\in I}\Spec\,\Z[Q_a].
$$
Hence, let us set
$$
U_a:=\Spec(\Z,Q_a)
\qquad
\text{for every $a\in I$}.
$$
We claim that these locally defined log structures glue to a
well defined log structure $\underline Q$ on the whole of
$S'_\tau$. Indeed, let $a,b\in I$; we have
$$
U_{a,b}:=\Spec\,\Z[Q_a]\cap\Spec\,\Z[Q_b]=\Spec\,\Z[Q_a\otimes_PQ_b]
$$
and it is easily seen that $Q_a\otimes_PQ_b=Q_a[b^{-1}a]$,
{\em i.e.} the localization of $Q_a$ obtained by
inverting its element $a^{-1}b$, and this is of course
the same as $Q_b[a^{-1}b]$. Then lemma \ref{lem_localize-const}
implies that the log structures of $U_a$ and $U_b$
agree on $U_{a,b}$, whence the contention. It is
easy to check that the resulting log scheme is precisely
$\Bl_{IP^{\log}_S}(S, P^{\log}_S)$ : the details shall be
left to the reader.
\end{example}

\begin{proposition}\label{prop_existenz}
Let $(X,\underline M)$ be a log scheme with quasi-coherent log
structure, and $\cI\subset\underline M$ a coherent ideal. Then,
the morphism \eqref{eq_blow-up-morph} is a logarithmic blow up
of $\cI$.
\end{proposition}
\begin{proof} Let $f:(Y,\underline N)\to(X,\underline M)$
be a morphism of log schemes, and suppose that $\cI\underline N$
is an invertible ideal of $\underline N$; in this case, we have
already remarked (see \eqref{subsec_first-invert}) that the
projection
$\pi_{(Y,\underline N,\cI\underline N)}:
\Bl_{\cI\underline N}(Y,\underline N)\to(Y,\underline N)$
is an isomorphism. We deduce a morphism :
\set\begin{equation}\label{eq_was-psi}
\phi(f,\cI)\circ\pi_{(Y,\underline N,\cI\underline N)}^{-1}:
(Y,\underline N)\to\Bl_\cI(X,\underline M).
\end{equation}
To conclude, it remains to show that \eqref{eq_was-psi} is
the only morphism of log schemes whose composition with
$\pi_{(X,\underline M,\cI)}$ equals $f$. The latter assertion
can be checked locally on $X_\tau$, hence we may assume that
$\underline M$ admits a chart $P_X\to\underline M$, such that
$\cI=I\underline M$ for some finitely generated ideal $I\subset P$.
In this case, in view of \eqref{eq_special-isoblow}, the set
of morphisms of $(X,\underline M)$-schemes
$(Y,\underline N)\to\Bl_\cI(X,\underline M)$ is in natural
bijection with the set of $(S,P_S^{\log})$-morphisms
$(Y,\underline N)\to B:=\Bl_{IP^{\log}_S}(S, P^{\log}_S)$
(notation of example \ref{ex_explicit-blowup}). In other words,
we may assume that $(X,\underline M)=(S,P_S^{\log})$, and
$\cI=IP^{\log}_S$. Then $f$ is determined by $\log f$,
{\em i.e.} by a map $\beta:P\to\underline N(Y)$.
Let $a_1,\dots,a_k$ be a system of generators for $I$;
for each $i=1,\dots,k$, we let $U_i\subset Y$ be the subset
of all $y\in Y$ for which there exists a $\tau$-point
$\xi$ of $Y$ with $|\xi|=y$ and
\set\begin{equation}\label{eq_local-generator}
a_i\underline N{}_\xi=I\underline N{}_\xi.
\end{equation}
Notice that if $y\in U_i$, then \eqref{eq_local-generator}
holds for every $\tau$-point $\xi$ of $Y$ localized at $y$
(details left to the reader).

\begin{claim}\label{cl_a-bit-better}
The subset $U_i$ is open in $Y_i$ for every $i=1,\dots,k$,
and $Y=\bigcup_{i=1}^kU_i$.
\end{claim}
\begin{pfclaim} Say that $y\in U_i$, and let $\xi$ be a
$\tau$-point of $Y$ localized at $y$, such that
\eqref{eq_local-generator} holds. This means that, for
every $j=1,\dots,k$, there exists $u_j\in\underline N{}_\xi$
such that $a_j=u_ja_i$. Then, we may find a $\tau$-neighborhood
$h:U'\to $ of $\xi$ such that this identity persists in
$\underline N(U')$; thus, $h(U')\subset U_i$. Since $h$ is an
open map, this shows that $U_i$ is an open subset.

Next, let $\xi$ be any $\tau$-point of $Y$; by assumption,
we have $I\underline N{}_\xi=b\underline N{}_\xi$ for some
$b\in\underline N{}_\xi$; this means that for every $i=1,\dots,k$
there exists $u_i\in\underline N{}_\xi$ such that $a_i=u_ib$.
Since $a_1,\dots,a_k$ generate $I$, we must have
$u_i\in\underline N{}^\times_\xi$ for at least an index
$i\leq k$, in which case $|\xi|\in U_i$, and this shows that
the $U_i$ cover the whole of $Y$, as claimed.
\end{pfclaim}

It is easily seen that, for every $i=1,\dots,k$, any morphism
$(U_i,\underline N{}_{|U_i})\to B$ of $(S,P^{\log}_S)$-schemes
factors through the open immersion $\Spec\,(Z,Q_{a_i})\to B$
(where $Q_a$, for an element $a\in P$, is defined as in
example \ref{ex_explicit-blowup}).
Conversely, by construction $\beta$ extends to a unique morphism
of monoids $Q_{a_i}\to\underline N(U_i)$. Summing up, there exists
at most one morphism of $(S,P_S^{\log})$-schemes
$(U_i,\underline N{}_{|U_i})\to B$. In light of claim
\ref{cl_a-bit-better}, the proposition follows.
\end{proof}

\sset\subsubsection{}\label{subsec_sat-blowup}
Keep the notation of proposition \ref{prop_existenz}; by inspecting
the construction, it is easily seen that the log structure
$\underline M'$ of $\Bl_\cI(X,\underline M)$ is quasi-coherent, and
if $\underline M$ is coherent (resp. quasi-fine, resp. fine), then
the same holds for $\underline M'$.
However, simple examples show that $\underline M'$ may fail to be
saturated, even in cases where $(X,\underline M)$ is a fs log
scheme. Due to the prominent role played by fs log schemes, it is
convenient to introduce the special notation :
$$
\satBl_\cI(X,\underline M):=(\Bl_\cI(X,\underline M))^\mathrm{qfs}
$$
for the {\em saturated logarithmic blow up\/} of a coherent ideal
$\cI$ in a quasi-fine log structure $\underline M$ (notation of
proposition \ref{prop_quasi-fine}).
Clearly the projection $\satBl_\cI(X,\underline M)\to(X,\underline M)$
is a final object of the category of saturated $(X,\underline M)$-schemes
in which the preimage of $\cI$ is invertible. If $(X,\underline M)$
is a fine log scheme, $\satBl_\cI(X,\underline M)$ is a fs log scheme.
Moreover, for any morphism of schemes $f:Y\to X$, let
$\underline M_Y:=f^*\underline M$; from remarks \eqref{rem_quite}(ii)
and \ref{rem_quasi-fine}(iii), we deduce a natural
isomorphism of $(Y,\underline M_Y)$-schemes :
\set\begin{equation}\label{eq_quick}
\satBl_{\cI\underline M_Y}(Y,\underline M_Y)\isom
Y\times_X\satBl_\cI(X,\underline M).
\end{equation}

\begin{theorem}\label{th_elementary-blowup}
Let $(X,\underline M)$ be a quasi-fine log scheme with saturated log
structure, $\cI\subset\underline M$ an ideal, $\xi$ a $\tau$-point of
$X$. Suppose that, in a neighborhood of $\xi$, the ideal $\cI$ is
generated by at most two sections, and denote :
$$
\phi:\Bl_\cI(X,\underline M)\to(X,\underline M)
\qquad
\text{(resp. $\phi_\sat:\satBl_\cI(X,\underline M)\to(X,\underline M)$)}
$$
the logarithmic (resp. saturated logarithmic) blow up of $\cI$. Then :
\begin{enumerate}
\item
If $\cI_\xi$ is an invertible ideal of $\underline M_\xi$,
the natural morphisms
$$
\phi^{-1}(\xi)\to\Spec\,\kappa(\xi)
\qquad
\phi_\sat^{-1}(\xi)\to\Spec\,\kappa(\xi)
$$
are isomorphisms.
\item
Otherwise, $\phi^{-1}(\xi)$ is a $\kappa(\xi)$-scheme isomorphic
to $\P^1_{\kappa(\xi)}$; furthermore, the same holds for the
reduced fibre $\phi_\sat^{-1}(\xi)_\mathrm{red}$, provided
$\underline M$ is fine.
\end{enumerate}
\end{theorem}
\begin{proof} After replacing $X$ by a $\tau$-neighborhood of
$\xi$, we reduce to the case where $\underline M$ admits an
integral and saturated chart $\alpha:P_X\to\underline M$
(lemma \ref{lem_simple-charts-top}(iii)), and if $\underline M$
is a fs log structure, we may also assume that the chart
$\alpha$ is fine and sharp at $\xi$ (corollary
\ref{cor_sharpy}(i)). Furthermore, we may assume that $\cI$ is
generated by at most two elements of $\underline M(X)$, and
if $\cI_\xi$ is principal, we may assume that the same holds
for $\cI$. In the latter case, since $\underline M$ is integral,
$\cI$ is invertible, hence $\phi$ and $\phi_\sat$ are isomorphisms,
so (i) follows already.

Now, suppose that $\cI_\xi$ is not invertible, and let
$a',b'\in\underline M(X)$ be a system of generators for $\cI$;
we can write $a'=\alpha(a)\cdot u$, $b'=\alpha(b)\cdot v$ for
some $a,b\in P$ and $u,v\in\kappa^\times$. Set $t:=a^{-1}b$,
and let $J\subset P$ be the ideal generated by $a$ and $b$;
clearly $J\underline M=\cI$, and $Pa,Pb\neq J$.

Consider first the case where $X=\Spec\,\kappa$, where $\kappa$
is a field (resp. a separably closed field, in case $\tau=\et$).
In this situation, a pre-log structure on $X_\tau$ is the same
as a morphism of monoids $\beta:P\to\kappa$, the associated
log structure is the induced map of monoids
$$
\beta^{\log}:P\otimes_{P_0}\kappa^\times
\to\kappa
\qquad\text{where}\qquad
P_0:=\beta^{-1}\kappa^\times
$$
and $\alpha$ is the natural map $P\to P\otimes_{P_0}\kappa^\times$.
After replacing $P$ by its localization $P_0^{-1}P$, we may also
assume that $P_0=P^\times$. Let
$$
(S,P^{\log}_S):=\Spec(\kappa,P)
\qquad
J^\sim:=JP^{\log}_S
\qquad
(Y,\underline N):=\Bl_{J^\sim}(S,P^{\log}_S).
$$
Denote by $\eps:\kappa[P]\to\kappa$ the homomorphism of
$\kappa$-algebras induced by $\beta$ via the adjunction
\eqref{subsec_mon-to-algs}, and set $I:=\Ker\,\eps$.
In view of lemma \ref{lem_blow-alg}(i) and \eqref{eq_quick},
we have natural cartesian diagrams of $\kappa$-schemes :
\set\begin{equation}\label{eq_fibre-of-blow}
{\diagram
\phi^{-1}(\xi) \ar[rr] \ar[d] & & (Y,\underline N) \ar[d] &
\phi_\sat^{-1}(\xi) \ar[rr] \ar[d] & &
(Y,\underline N)^\sat \ar[d] \\
|\xi| \ar[rr]^-{\Spec\,\eps} & & S &
|\xi| \ar[rr]^-{\Spec\,\eps} & & S.
\enddiagram}
\end{equation}
On the other hand, $(a,b)$ can be regarded as a pair of global
sections of $J^\sim\underline N$, so the universal property of
example \ref{ex_proj-space-continued}(ii) yields a morphism of
$(S,P^{\log}_S)$-schemes :
$$
f_{(a,b)}:(Y,\underline N)\to\P^1_{(S,P^{\log}_S)}.
$$
In light of example \ref{ex_proj-space-continued}(i), the
assertion concerning $\phi^{-1}(\xi)$ will then follow from the :

\begin{claim}\label{cl_excluded}
The morphism $|\xi|\times_Sf_{(a,b)}$ is an isomorphism
of $\kappa$-schemes.
\end{claim}
\begin{pfclaim} Let $Q_a\subset P^\gp$ (resp. $Q_b\subset P^\gp$)
be the submonoid generated by $P$ and $t$ (resp. by $P$ and
$t^{-1}$); by inspecting example \ref{ex_explicit-blowup}, we see
that $Y$ is covered by two affine open subsets :
$$
U_a:=\Spec\,\kappa[Q_a]
\qquad
U_b:=\Spec\,\kappa[Q_b]
$$
and $U_a\cap U_b=\Spec\,\kappa[Q_a\otimes_PQ_b]$.
On the other hand, $\P^1_{(S,P^{\log}_S)}$ is covered as
well by two affine open subsets $U'_0$ and $U'_\infty$, both
isomorphic to $\Spec\,\kappa[P\oplus\N]$, and such that
$U'_0\cap U'_\infty=\Spec\,\kappa[P\oplus\Z]$, as usual.
Moreover, a direct inspection shows that $f_{(a,b)}$ restricts
to morphisms
$$
U_a\to U'_0
\qquad
U_b\to U'_\infty
$$
induced respectively by the maps of $\kappa$-algebras
$$
\omega_0:\kappa[P\oplus\N]\to\kappa[Q_a]
\qquad
\omega_\infty:\kappa[P\oplus\N]\to\kappa[Q_b]
$$
such that $\omega_0(x,n)=x\cdot t^n$ and
$\omega_\infty(x,n)=x\cdot t^{-n}$ for every $(x,n)\in P\oplus\N$.
We show that $\bar\omega_0:=\omega_0\otimes_{\kappa[P]}\kappa[P]/I$
is an isomorphism; the same argument will apply also to
$\omega_\infty$, so the claim shall follow.

Indeed, clearly the $\kappa[P]$-algebra $\kappa[Q_a]$
is generated by $t$, hence $\omega_0$ is surjective, and then the
same holds for $\bar\omega_0$. Next, set $H_a:=I\cdot\kappa[Q_a]$;
it is easily seen that $H_a$ consists of all sums of the
form $\sum_{j=0}^nc_jt^j$, for arbitrary $n\in\N$, with
$c_j\in I$ for every $j=0,\dots,n$. Clearly, an element
$p(T)\in\kappa[\N]=\kappa[T]$ lies in $\Ker\,\bar\omega_0$
if and only if $p(t)\in H_a$, so we come down to the following
assertion. Let $c_0,\dots,c_n\in\kappa[P]$ such that
\set\begin{equation}\label{eq_little-sum}
\sum_{j=0}^nc_jt^j=0
\end{equation}
in $\kappa[Q_a]$; then $c_j$ lies in the ideal
$\kappa[\beta^{-1}(0)]$ of $\kappa[P]$, for every $j=0,\dots,n$.
Since $P$ is integral, \eqref{eq_little-sum} is equivalent
to the identity : $\sum_{j=0}^nc_ja^{n-j}b^j=0$ in $\kappa[P]$.
For every $x\in P$, denote by $\pi_x:\kappa[P]\to\kappa$
the $\kappa$-linear projection such that $\pi_x(x)=1$, and
$\pi_x(y)=0$ for every $y\in P\setminus\{x\}$.

Suppose, by way of contradiction, that
$c_i\notin\kappa[\beta^{-1}(0)]$ for some $i\leq n$, hence
$\pi_x(c_i)\neq 0$ for some $x\in P_0$; since $P_0=P^\times$,
we may replace $c_j$ by $x^{-1}c_j$, for every $j\leq n$, and
assume that $\pi_1(c_i)\neq 0$, hence
$\pi_{a^{n-i}b^i}(c_ia^{n-i}b^i)\neq 0$ (again, using the
assumption that $P$ is integral). Thus, there exists $j\neq i$
with $j\leq n$, such that $\pi_{a^{n-i}b^i}(c_ja^{n-j}b^j)\neq 0$,
and we may then find an element $c\in P$ such that
$\pi_{a^{n-i}b^i}(ca^{n-j}b^j)=1$, {\em i.e.}
$ca^{n-j}b^j=a^{n-i}b^i$; up to swapping the roles of $a$
and $b$, we may assume that $i>j$, in which case we may
write $t^{i-j}=c$; since $P$ is saturated, it follows that
$t\in P$, hence $J$ is generated by $a$, which is excluded.
\end{pfclaim}

Next, we assume that $\underline M$ is a fs log structure
(and $X$ is still $\Spec\,\kappa$), and we consider the
morphism $\phi_\sat$. As already remarked, we may assume
that $\alpha$ is sharp at $\xi$, and $P$ fine and saturated;
the sharpness condition amounts to saying that $\beta(x)=0$ for
every $x\in P\setminus\{1\}$, therefore $I$ is the augmentation
ideal of the graded $\kappa$-algebra $\kappa[P]$. By inspecting
the proof of proposition \ref{prop_quasi-fine}, we see that
$(Y,\underline N)^\sat$ is covered by two affine open subsets
$$
U_a^\mathrm{fs}:=\Spec\,\kappa[Q_a^\sat]
\qquad
U_b^\mathrm{fs}:=\Spec\,\kappa[Q_b^\sat]
$$
and $U_a^\mathrm{fs}\cap U_b^\mathrm{fs}=
\Spec\,\kappa[Q_a^\sat\otimes_PQ_b^\sat]$.
Since $J$ is not principal, we have $t\notin P$, and since $P$
is saturated, we deduce that $t$ is not a torsion element
of $P^\gp$; as the latter is a finitely generated abelian
group, it follows that we may find a {\em unimodular}
element $u\in P^\gp$ such that $t$ lies in the submonoid
$\N u\subset P^\gp$ generated by $u$; this condition means
that $t=u^k$ for some $k\in\N$, and $\N u$ is not properly
contained in another rank one free submonoid of $P^\gp$.
Write $u=a'{}^{-1}b'$ for some $a',b'\in P$, let $J'\subset P$
be the ideal generated by $a'$ and $b'$, and $R_{a'}$ (resp.
$R_{b'}$) the submonoid of $P^\gp$ generated by $P$ and $u$
(resp. by $P$ and $u^{-1}$); clearly $R_{a'}^\sat=Q_a^\sat$,
and $R_{b'}^\sat=Q_b^\sat$. Denote by $\underline N'$ the log
structure of $(Y,\underline N)^\sat$, and $\cJ':=J'\underline N'$;
it is easily seen that
$$
\cJ'_{|U^\mathrm{fs}_a}=a'\underline N'_{|U^\mathrm{fs}_a}
\qquad
\cJ'_{|U^\mathrm{fs}_b}=b'\underline N'_{|U^\mathrm{fs}_b}
$$
hence $\cJ'$ is invertible, and example
\ref{ex_proj-space-continued}(ii) yields a morphism of
$(S,P^{\log}_S)$-schemes
$$
f_{(a',b')}:(Y,\underline N)^\sat\to\P^1_{(S,P^{\log}_S)}.
$$
In light of example \ref{ex_proj-space-continued}(i), the
assertion concerning $\phi^{-1}_\sat(\xi)$ will then follow
from the :

\begin{claim}\label{cl_comm-avant}
$(|\xi|\times_Sf_{(a',b')})_\red:
\phi^{-1}_\sat(\xi)_\red\to\P^1_\kappa$ is an isomorphism of
$\kappa$-schemes.
\end{claim}
\begin{pfclaim} As in the proof of claim \ref{cl_excluded},
the morphism $f_{(a',b')}$ restricts to morphisms
$U^\mathrm{fs}_a\to U'_0$ and $U^\mathrm{fs}_b\to U'_\infty$
induced by maps of $\kappa$-algebras :
$$
\omega_0:\kappa[P\oplus\N]\to\kappa[Q_a^\sat]
\qquad
\omega_\infty:\kappa[P\oplus\N]\to\kappa[Q_b^\sat]
$$
such that $\omega_0(x,n)=x\cdot u^n$ and
$\omega_\infty(x,n)=x\cdot u^{-n}$ for every $(x,n)\in P\oplus\N$,
and again, it suffices to show that
$$
(\omega_0\otimes_{\kappa[P]}\kappa[P]/I)_\red:
\kappa[T]\to(\kappa[R^\sat_a]/I\kappa[R^\sat_a])_\red
$$
is an isomorphism (where for any ring $A$, we denote
$A_\mathrm{red}$ be the maximal reduced quotient of $A$,
{\em i.e.} $A_\mathrm{red}:=A/\nil(A)$, where $\nil(A)$
is the nilradical of $A$). We break the latter
verification in two steps : first, let us check that
the map
$$
\bar\omega'_0:\kappa[T]\to\kappa[R_a]/I\kappa[R_a]
\qquad
p(T)\mapsto p(u) \pmod{I\kappa[R_a]}
$$
is an isomorphism. Indeed, $\bar\omega'_0$ is induced by the map
of monoids $\phi:\N\to R_a$ such that $n\mapsto u^n$ for every
$n\in\N$; if $t=u^k$, the map $\phi$ fits into the cocartesian
diagram :
$$
\xymatrix{ \N \ar[r]^-\psi \ar[d]_{\bek_\N} & Q_a \ar[d]^j \\
           \N \ar[r]^-\phi & R_a
}$$
where $\bek_\N$ is the $k$-Frobenius map, $\psi$ is given by
the rule : $n\mapsto t^n$ for every $n\in\N$, and $j$ is
the natural inclusion. Hence :
$$
\bar\omega'_0=(\kappa[\psi]\otimes_{\kappa[P]}\kappa[P]/I)
\otimes_{\kappa[T^k]}\kappa[T].
$$
However, the proof of claim \ref{cl_excluded} shows that
$\kappa[\psi]\otimes_{\kappa[P]}\kappa[P]/I$ is an isomorphism,
whence the contention. Lastly, let show that the natural map
$$
\bar\omega''_0:\kappa[R_a]/I\kappa[R_a]\to
(\kappa[R_a^\sat]/I\kappa[R_a^\sat])_\mathrm{red}
$$
is an isomorphism. Indeed, it is clear that the natural map
$\omega''_0:\kappa[R_a]\to\kappa[R_a^\sat]$ is integral and
injective, hence $\Spec\,\omega''_0$ is surjective; therefore
$\Spec\,\bar\omega''_0$ is still surjective and integral.
However, the foregoing shows that $\kappa[R_a]/I\kappa[R_a]$
is reduced, so we deduce that $\bar\omega''_0$ is injective.
To show that $\bar\omega''_0$ is surjective, it suffices to
show that the classes of the generating system
$R^\sat_a\subset\kappa[R_a^\sat]$ lie in the image of $\bar\omega''_0$.
Hence, let $x\in R_a^\gp$, with $x^m\in R_a$ for some $m>0$, so that
$x^m=y\cdot u^n$ for some $n\in\N$ and $y\in P$. If $y\neq 1$, we
have $y\in I$, hence the image of $x^m$ vanishes in
$\kappa[R_a^\sat]/I\kappa[R_a^\sat]$, and the image of $x$ vanishes
in the reduced quotient; finally, if $y=1$, the identity $x^m=u^n$
implies that $m$ divides $n$, since $u$ is unimodular; hence
$x=u^{n/m}$ and the image of $x$ agrees with $\bar\omega''_0(u^{n/m})$.
\end{pfclaim}

Finally, let us return to a general quasi-fine log scheme
$(X,\underline M)$; the theorem will follow from the more precise :

\begin{claim}\label{cl_more-and-more}
In the situation of the theorem, suppose moreover that
\begin{enumerate}
\alphaenu
\item
$\underline M$ admits a saturated chart $\alpha:P_X\to\underline M$
\item
$\cI=J\underline M$, where $J\subset P$ is an ideal generated by two
elements $a,b\in P$
\item
$\cI_\xi$ is not invertible.
\end{enumerate}
Then we have :
\begin{enumerate}
\item
There exists a morphism of $(X,\underline M)$-schemes :
$\Bl_\cI(X,\underline M)\to\P^1_{(X,\underline M)}$\ \
inducing an isomorphism of $\kappa(\xi)$-schemes
$\phi^{-1}(\xi)\isom\P^1_{\kappa(\xi)}$.
\item
If furthermore, $P$ is fine (and saturated) and $\alpha$
is sharp at $\xi$, then there exists a morphism of
$(X,\underline M)$-schemes :
$\satBl_\cI(X,\underline M)\to\P^1_{(X,\underline M)}$\ \
inducing an isomorphism of $\kappa(\xi)$-schemes
$\phi^{-1}_\sat(\xi)_\red\isom\P^1_{\kappa(\xi)}$.
\end{enumerate}
\end{claim}
\begin{pfclaim}[](i): Denote by $\underline N$ the log structure of
$\Bl_\cI(X,\underline M)$; the elements $a,b$ define global
sections of the invertible $\underline N$-module $\cI\underline N$,
and we claim that the corresponding morphism of
$(X,\underline M)$-schemes
$f_{(a,b)}:\Bl_\cI(X,\underline M)\to\P^1_{(X,\underline M)}$
will do. Indeed, set
$(|\xi|,\underline M{}_\xi):=|\xi|\times_X(X,\underline M)$,
and recall that there exists natural isomorphisms
$$
|\xi|\times_X\P^1_{(X,\underline M)}\isom
\P^1_{(|\xi|,\underline M{}_\xi)}
\qquad
|\xi|\times_X\Bl_\cI(X,\underline M)\isom
\Bl_{\cI\underline M{}_\xi}(|\xi|,\underline M{}_\xi)
$$
(example \ref{ex_proj-space-continued}(i) and remark
\ref{rem_quite}(ii)). Denote by $\underline N{}_\xi$ the
log structure of
$\Bl_{\cI\underline M{}_\xi}(|\xi|,\underline M{}_\xi)$;
by example \ref{ex_proj-space-continued}(iii), the base change
$$
|\xi|\times_Xf:\Bl_{\cI\underline M{}_\xi}(|\xi|,\underline M{}_\xi)
\to\P^1_{(|\xi|,\underline M{}_\xi)}
$$
is the unique morphism $f_{(\bar a,\bar b)}$ of
$(|\xi|,\underline M{}_\xi)$-schemes corresponding to the pair
$(\bar a,\bar b)$ of global sections of $\cI\underline N{}_\xi$
obtained by pulling back the pair $(a,b)$. Therefore, in order
to check that $|\xi|\times_Xf$ is an isomorphism, we may replace
from start $(X,\underline M)$ by $(|\xi|,\underline M{}_\xi)$
(whose log structure is still quasi-fine, by lemma
\ref{lem_simple-charts-top}(i), and assume that $X=\Spec\,\kappa$,
where $\kappa$ is a field (resp. a separably closed field, in case
$\tau=\et$), in which case the assertion is just claim
\ref{cl_excluded}.

(ii): Denote by $\underline N'$ the log structure of
$\satBl_\cI(X,\underline M)$, define $a',b'$ and $J'$ as
in the foregoing, and set again $\cJ':=J'\underline N'$.
Again, it is easily seen that $\cJ'$ is an invertible
$\underline N'$-module, and the pair $(a',b')$ yields
a morphism
$f_{(a',b')}:\satBl_\cI(X,\underline M)\to\P^1_{(X,\underline M)}$
which fulfills the sought condition. Indeed, denote
by $\underline N{}'_\xi$ the log structure of
$\satBl_{\cI\underline M{}_\xi}(|\xi|,\underline M{}_\xi)$;
in light of \eqref{eq_quick} and example
\ref{ex_proj-space-continued}(iii), the base change
$$
|\xi|\times_Xf_{(a',b')}:
\satBl_{\cI\underline M{}_\xi}(|\xi|,\underline M{}_\xi)
\to\P^1_{(|\xi|,\underline M{}_\xi)}
$$
is the unique morphism $f_{(\bar a{}',\bar b{}')}$ of
$(|\xi|,\underline M{}_\xi)$-schemes corresponding to the pair
$(\bar a{}',\bar b{}')$ of global sections of
$\cJ'\underline N{}'_\xi$ obtained by pulling back the pair
$(a',b')$. Thus, the assertion is just claim \ref{cl_comm-avant}.
\end{pfclaim}
\end{proof}

\subsection{Regular log schemes}\label{sec_log-regular}
In this section we introduce the logarithmic version of the
classical regularity condition for locally noetherian schemes.
This theory is essentially due to K.Kato (\cite{Ka2}), and we
mainly follow his exposition, except in a few places where his
original arguments are slightly flawed, in which cases we supply
the necessary corrections.

\sset\subsubsection{}\label{subsec_Krull-dim-mon}
Let $A$ be a ring, $P$ a monoid. Recall that $\fm_P$ is the maximal
(prime) ideal of $P$ (see \eqref{subsec_sepc-of-monoid}). The
{\em $\fm_P$-adic filtration\/} of $P$ is the descending sequence
of ideals :
$$
\dots\subset\fm_P^3\subset\fm_P^2\subset\fm_P
$$
where $\fm^n_P$ is the $n$-th power of $\fm$ in the monoid $\cP(P)$
(see \eqref{subsec_toric}). It induces a {\em $\fm_P$-adic
filtration\/} $\Fil_\bullet M$ on any $P$-module $M$ and any
$A[P]$-algebra $B$, defined by letting $\Fil_nM:=\fm^n_PM$ and
$\Fil_nB:=A[\fm^n_P]\cdot B$, for every $n\in\N$.

\begin{lemma}\label{lem_canonickal}
Suppose that $P$ is fine. We have :
\begin{enumerate}
\item
The $\fm_P$-adic filtration is separated on $P$.
\item
If $P$ is sharp, $P\!\setminus\!\fm^n_P$ is a finite set, for every
$n\in\N$.
\end{enumerate}
\end{lemma}
\begin{proof} (i): Indeed, choose $A$ to be a non-zero noetherian
ring, set $J:=\bigcap_{n\geq 0}A[\fm^n_P]$ and notice that $J$ is
generated by $\fm^\infty_P:=\bigcap_{n\in\N}\fm^n_P$. On the other
hand, $J$ is annihilated by an element of $1+A[\fm_P]$
(\cite[Th.8.9]{Mat}).
Thus, suppose $x\in\fm^\infty_P$, and pick $y\in A[\fm_P]$ such that
$(1-y)x=0$; we may write $y=a_1t_1+\cdots+a_rt_r$ for certain
$a_1,\dots,a_r\in A$ and $t_1,\dots,t_r\in\fm_P$. Therefore
$x=a_1xt_1+\cdots+a_rxt_r$ in $A[P]$, which is absurd, since $P$ is
integral.

Assertion (ii) is immediate from the definition.
\end{proof}

\sset\subsubsection{}\label{subsec_noncanonick}
Keep the assumptions of lemma \ref{lem_canonickal}. It turns out
that $P$ can actually be made into a graded monoid, albeit in a
non-canonical manner.
We proceed as follows. Let $\eps:P\to P^\sat$ the inclusion map,
$T\subset P^\sat$ the torsion subgroup, set $Q:=P^\sat/T$, and let
$\pi:P^\sat\to Q$ be the natural surjection. We may regard $\log Q$
as a submonoid of the polyhedral cone $Q_\R$, lying in the vector
space $Q^\gp_\R$, as in \eqref{subsec_Gordon-more}. Since $Q_\R$ is
a rational polyhedral cone, the same holds for $Q^\vee_\R$, hence we
may find a $\Q$-linear form $\gamma:\log Q^\gp\otimes_\Z\Q\to\Q$,
which is non-negative on $\log Q$, and such that
$Q_\R\cap\Ker\,\gamma\otimes_\Q\R$ is the minimal face of $Q_\R$, {\em
i.e.} the $\R$-vector space spanned by the image of $Q^\times$. If
we multiply $\gamma$ by some large positive integer, we may achieve
that $\gamma(\log P)\subset\N$. We set :
$$
\gr^\gamma_n P:=(\gamma\circ\pi\circ\eps)^{-1}(n) \qquad \text{for
every $n\in\N$.}
$$
It is clear that
$\gr^\gamma_nP\cdot\gr^\gamma_mP\subset\gr^\gamma_{n+m}P$, hence
$$
P=\coprod_{n\in\N}\gr_n^\gamma P
$$
is a $\N$-graded monoid, and consequently :
$$
A[P]=\bigoplus_{n\in\N} A[\gr^\gamma_nP]
$$
is a graded $A$-algebra. Moreover, it is easily seen that
$\gr_0^\gamma P=P^\times$. More generally, for $x\in P$,
let $\mu(x)$ be the maximal $n\in\N$ such that $x\in\fm^n_P$;
then there exists a constant $C\geq 1$ such that :
$$
\gamma(x)\geq\mu(x)\geq C^{-1}\gamma(x) \qquad \text{for every $x\in
P$}
$$
so that the $\fm_P$-adic filtration and the filtration defined by
$\gr^\gamma_\bullet P$, induce the same topology on $P$ and on $A[P]$.
As a corollary of these considerations, we may state the following
``regularity criterion'' for fine monoids :

\begin{proposition}\label{prop_regul-critt}
Let $P$ be an integral monoid such that $P^\sharp$ is fine.
Then we have
$$
\rk^\circ_{P^\times}\fm_P/\fm_P^2\geq\dim P
$$
(notation of example {\em\ref{ex_rank-free-point-mod}}) and the
equality holds if and only if $P^\sharp$ is a free monoid.
\end{proposition}
\begin{proof} (Notice that $\fm_P/\fm_P^2$ is a free pointed
$P^\times$-module, since $P^\times$ obviously acts freely
on $\fm_P\setminus\fm_P^2$.) Since
$\fm_{P^\sharp}/\fm_{P^\sharp}^2=(\fm_P/\fm_P^2)\otimes_PP^\sharp$,
we may replace $P$ by $P^\sharp$, and assume from start that
$P$ is sharp and fine. Then, the rank of $\fm_P/\fm_P^2$ equals
the cardinality of the set $\Sigma:=\fm_P\setminus\fm_P^2$, which
is finite, by lemma \ref{lem_canonickal}. We have a surjective
morphism of monoids $\phi:\N^{(\Sigma)}\to P$, that sends the
basis of $\N^{(\Sigma)}$ bijectively onto $\Sigma\subset P$
(corollary \ref{cor_ideals-in-fg-mon}).
The sought inequality follows immediately, and it is also clear
that we have equality, in case $P$ is free. Conversely, if
equality holds, $\phi^\gp$ must be a surjective group homomorphism
between free abelian group of the same finite rank (corollary
\ref{cor_consequent}(i)), so it is an isomorphism, and then the
same holds for $\phi$.
\end{proof}

\begin{proposition}\label{prop_Krull-dim-mon}
Let $P$ be a fine and sharp monoid, $A$ a noetherian local ring.
Set $S_P:=1+A[\fm_P]$; then we have :
$$
\dim S_P^{-1}A[P]=\dim A+\dim P.
$$
\end{proposition}
\begin{proof} To begin with, the assumption that $P$ is sharp
implies that $S_P^{-1}A[P]$ is local. Next, notice that $A[P]$
is a free $A$-module, hence the natural map $A\to S_P^{-1}A[P]$
is a flat and local ring homomorphism. Let $k$ be the residue field
of $A$; in view of \cite[Th.15.1(ii)]{Mat}, we deduce :
$$
\dim S_P^{-1}A[P]=\dim A+\dim S_P^{-1}k[P].
$$
Hence we are reduced to showing the stated identity for $A=k$. 
However, clearly $S_P^{-1}k[P]=k[P]_\fm$, where $\fm$ is the maximal
ideal generated by the image of $\fm_P$, hence it suffices to apply
claim \ref{cl_NullSt}(ii) and corollary \ref{cor_consequent}(i), to
conclude.
\end{proof}

\sset\subsubsection{}\label{subsec_constant-term}
Let $A$ be a ring, and $P$ a fine and sharp monoid. We define :
$$
A[[P]]:=\lim_{n\in\N}\, A\La P/\fm^n_P\Ra.
$$
Alternatively, this is the completion of $A[P]$ for its
$A[\fm_P]$-adic topology.
In view of the finiteness properties of the $\fm_P$-adic filtration
(lemma \ref{lem_canonickal}(ii)), one may present $A[[P]]$ as the
ring of formal infinite sums $\sum_{\sigma\in
P}a_\sigma\cdot\sigma$, with arbitrary coefficients $a_\sigma\in A$,
where the multiplication and addition are defined in the obvious
way. Moreover, we may use a morphism of monoids $\gamma:\log P\to\N$
as in \eqref{subsec_noncanonick}, to see that :
\set\begin{equation}\label{eq_decomponi}
A[[P]]=\prod_{n\in\N}A[\gr^\gamma_nP]
\end{equation}
where $A[\gr^\gamma_nP]\cdot A[\gr^\gamma_mP]\subset
A[\gr^\gamma_{n+m}P]$ for every $m,n\in\N$. So any element $x\in
A[[P]]$ can be decomposed as an infinite sum
$$
x=\sum_{n\in\N}\gr^\gamma_n x.
$$
The term $\gr^\gamma_0 x\in\gr_0^\gamma A=A$ does not depend on the
chosen $\gamma$ : it is the {\em constant term\/} of $x$, {\em i.e.}
the image of $x$ under the natural projection $A[[P]]\to A$.

\begin{corollary}\label{cor_new-dim}
Let $P$ be a fine and sharp monoid, $A$ a noetherian local ring.
Then :
\begin{enumerate}
\item
For any local morphism $P\to A$ (see \eqref{subsec_sepc-of-monoid}),
we have the inequality:
$$
\dim A\leq\dim A/\fm_PA+\dim P.
$$
\item
\ $\dim A[P]=\dim A[[P]]=\dim A+\dim P$.
\end{enumerate}
\end{corollary}
\begin{proof}(i): Set $A_0:=A/\fm_PA$, and $B:=S_P^{-1}A_0[P]$, where
$S_P\subset A_0[P]$ is the multiplicative subset $1+A_0[\fm_P]$; if
we denote by $\gr_\bullet A$ (resp. $\gr_\bullet B$) the graded
$A_0$-algebra associated to the $\fm_P$-adic filtration on $A$
(resp. on $B$), we have a natural surjective homomorphism of graded
$A_0$-algebras :
$$
\gr_\bullet B\to\gr_\bullet A.
$$
Hence $\dim A=\dim\gr_\bullet A\leq\dim\gr_\bullet B=\dim B$, by
\cite[Th.15.7]{Mat}. Then the assertion follows from proposition
\ref{prop_Krull-dim-mon}.

(ii): Set $B:=S_P^{-1}A[P]$ and let $\fm_A$ (resp $\fm_B$)
be the maximal ideal of $A$ (resp. of $B$); notice that
$A[[P]]$ is the $(\fm_PB)$-adic completion of the local
noetherian ring $B$, so $A[[P]]$ is a local noetherian
ring as well, with maximal ideal
$\fn:=\fm_A[[P]]+A[[\fm_P]]$, and the $\fn$-adic completion
of $A[[P]]$ is naturally isomorphic to the $\fm_B$-adic
completion $B^\wedge$ of $B$. Hence we get
$$
\dim A[[P]]=\dim B^\wedge=\dim B
$$
so the second identity follows from proposition
\ref{prop_Krull-dim-mon}. Next, clearly we have
$\dim A\geq\dim B$. On the other hand, let
$\fq\subset A[P]$ be any prime ideal, set $\fp:=\fq\cap A$
and denote by $\kappa$ the residue field of the
local ring $A_\fp$; with this notation,
\cite[Th.15.1(ii)]{Mat} says that
$$
\dim A[P]_\fq=\dim A_\fp+\dim A[P]_\fq\otimes_A\kappa
=\dim A_\fp+\dim\kappa[P]_\fq\leq\dim A+\dim\kappa[P]
$$
so the assertion follows from \ref{cl_NullSt}(ii) and
corollary \ref{cor_consequent}(i).
\end{proof}

As a first application, we have the following combinatorial
version of Kunz's theorem \ref{th_Kunz-by-Matsu} that characterizes
regular rings via their Frobenius endomorphism.

\begin{theorem}\label{th_Kunz}
Let $P$ be a monoid such that $P^\sharp$ is fine, $k>1$ an integer,
and suppose that the Frobenius endomorphism $\bek_P:P\to P$ is flat
(see example {\em\ref{ex_multiply-by-n-in-fan}(i)}). Then $P^\sharp$
is a free monoid.
\end{theorem}
\begin{proof} First, we remark that
$\bek^\sharp_P:P^\sharp\to P^\sharp$ is still flat (corollary
\ref{cor_yet-another-flat}(i)). Moreover, $\bek^\sharp_P$
is injective. Indeed, suppose that $x^k=y^k\cdot u$ for
some $x,y\in P$ and $u\in P^\times$; from theorem
\ref{th_always-exact} we deduce that there exist
$b_1,b_2,t\in P$ such that
$$
b_1x=b_2y
\qquad
1=b_1^kt
\qquad
u=b_2^kt.
$$
Especially, $b_1,b_2\in P^\times$, so the images of $x$ and
$y$ agree in $P^\sharp$. Hence, we may replace $P$ by $P^\sharp$,
and assume that $\bek_P$ is flat and injective, in which
case $\Z[\bek_P]:\Z[P]\to\Z[P]$ is flat (theorem
\ref{th_flat-crit-for-mnds}), integral and injective, hence
it is faithfully flat. Now, let $R$ be the colimit of
the system of rings $(R_n~|~n\in\N)$, where $R_n:=\Z[P]$, and
the transition map $R_n\to R_{n+1}$ is $\Z[\bek_P]$ for every
$n\in\N$. The induced map $j:R_0\to R$ is still faithfully flat;
moreover, let $p$ be any prime divisor of $k$, and notice that
$j\circ\bep_P=j$ (where $\bep_P$ is the $p$-Frobenius map). It
follows that $\bep_P$ is flat and injective as well, so
$\F_p[\bep_P]:\F_p[P]\to\F_p[P]$ is a flat ring homomorphism
(again, by theorem \ref{th_flat-crit-for-mnds}), and then the
same holds for the induced map $\F_p[[\bep_P]]:\F_p[[P]]\to\F_p[[P]]$.
By Kunz's theorem, we deduce that $\F_p[[P]]$ is a regular
local ring, with maximal ideal $\fm:=\F_p[[\fm_P]]$; notice
that the images of the elements of $\fm_P\!\setminus\!\fm_P^2$
yield a basis for the $\F_p$-vector space $\fm/\fm^2$.
Say that $\fm_P\!\setminus\!\fm_P^2=\{x_1,\dots,x_s\}$;
it follows that the continuous ring homomorphism
$$
\F_p[[T_1,\dots,T_s]]\to\F_p[[P]]
\qquad
T_i\mapsto x_i
\qquad
\text{for $i=1,\dots,s$}
$$
is an isomorphism. From the discussion of
\eqref{subsec_constant-term}, we immediately deduce that
$P\simeq\N^{\oplus s}$, as required.
\end{proof}

\sset\subsubsection{}\label{subsec_log-local-flat-crit}
Now we wish to state and prove the combinatorial versions of the
Artin-Rees lemma, and of the so-called local flatness criterion
(see {\em e.g.} \cite[Th.22.3]{Mat}).
Namely, let $P$ be a {\em pointed\/} monoid, such that
$P^\sharp$ is finitely generated; let also $(A,\fm_A)$ be a
local noetherian ring, $N$ a finitely generated $A$-module, and :
$$
\alpha:P\to(A,\cdot)
$$
a morphism of pointed monoids. The following is our version of the
Artin-Rees lemma :

\begin{lemma}\label{lem_log-Artin-Rees}
In the situation of \eqref{subsec_log-local-flat-crit}, let
$J\subset P$ be an ideal, $M$ a finitely generated $P$-module,
$M_0\subset M$ a submodule. Then there exists $c\in\N$ such that :
\set\begin{equation}\label{eq_log_Artin-Rees}
J^nM\cap M_0=J^{n-c}(J^cM\cap M_0) \qquad \text{for every $n>c$}.
\end{equation}
\end{lemma}
\begin{proof} Set $\bar M:=M/P^\times$, $\bar M_0 :=M_0/P^\times$
and $\bar J:=J/P^\times$, the set-theoretic quotients for the
respective natural $P^\times$-actions. Notice that
$\bar J$ is an ideal of $P^\sharp$ and $\bar M_0\subset\bar M$ is an
inclusion of $\bar P$-modules. Moreover, any set of generators of
the ideal $\bar J$ (resp. of the $P^\sharp$-module $\bar M_0$) lifts
to a set of generators for $J$ (resp. for the $P$-module $M_0$).
Furthermore, it is easily seen that \eqref{eq_log_Artin-Rees} is
equivalent to the identity $\bar J{}^n\bar M\cap\bar M_0=\bar
J{}^{n-c}(\bar J{}^c\bar M\cap\bar M_0)$. Hence, we are reduced to
the case where $P=P^\sharp$ is a finitely generated monoid. Then
$\Z[P]$ is noetherian, $\Z[M]$ is a $\Z[P]$-module of finite type,
and we notice that :
$$
\Z[J^nM\cap M_0]=J^n\Z[M]\cap\Z[M_0]
\qquad
\Z[J^{n-c}(J^cM\cap M_0)]=J^{n-c}(J^c\Z[M]\cap\Z[M_0]).
$$
Thus, the assertion follows from the standard Artin-Rees lemma
\cite[Th.8.5]{Mat}.
\end{proof}

\begin{proposition}\label{prop_flat-vrit}
In the situation of \eqref{subsec_log-local-flat-crit},
suppose moreover that $P^\sharp$ is fine (see remark
{\em\ref{rem_apparent-reasons}(vi)}), and let
$\fm_\alpha:=\alpha^{-1}\fm_A$. Then the following
conditions are equivalent :
\begin{enumerate}
\alphaenu
\item
$N$ is $\alpha$-flat (see definition {\em\ref{def_M-flatness}}).
\item
$\Tor_i^{\Z\La P\Ra}(N,\Z\La M\Ra)=0$ for every $i>0$ and every
integral pointed $P$-module $M$.
\item
$\Tor_1^{\Z\La P\Ra}(N,\Z\La P/\fm_\alpha\Ra)=0$.
\item
The natural map :
$$
(\fm^n_\alpha/\fm^{n+1}_\alpha)\otimes_PN\to\fm^n_\alpha
N/\fm^{n+1}_\alpha N
$$
is an isomorphism of $A$-modules, for every $n\in\N$ (notation
of \eqref{subsec_restrict-scalar-mon}).
\end{enumerate}
\end{proposition}
\begin{proof} The assertion (a)$\Leftrightarrow$(b) is just
a restatement of proposition \ref{prop_fla-criterion-point}(i),
and holds in greater generality, without any assumption on
either $A$ or the pointed integral monoid $P$.

As for the remaining assertions, let $S:=\alpha^{-1}(A^\times)$;
since the localization $P\to S^{-1}P$ is flat, the natural maps :
$$
\Tor^{\Z\La P\Ra}_i(N,\Z\La M\Ra)\to\Tor_i^{\Z\La
S^{-1}P\Ra}(N,\Z\La S^{-1}M\Ra)
$$
are isomorphisms, for every $i\in\N$ and every $P$-module $M$. Also,
notice that the two $P$-modules appearing in (d) are actually
$S^{-1}P$-modules (and the natural map is $\Z\La
S^{-1}P\Ra$-linear). Hence, we can replace everywhere $P$ by
$S^{-1}P$, which allows to assume that $\alpha$ is local, {\em i.e.}
$\fm_\alpha=\fm_P$.

Next, obviously (b)$\Rightarrow$(c).

(c)$\Rightarrow$(d): For every $n\in\N$, we have a short exact
sequence of pointed $P$-modules :
$$
0\to\fm^n_P/\fm^{n+1}_P\to P/\fm^{n+1}_P\to P/\fm^n_P\to 0.
$$
It is easily seen that $\fm^n_P/\fm^{n+1}_P$ is a free
$P/\fm_P$-module (in the category of pointed modules), so the
assumption implies that $\Tor_1^{\Z\La
P\Ra}(N,\fm^n_P/\fm^{n+1}_P)=0$ for every $n\in\N$. By looking at
the induced long $\Tor$-sequences, we deduce that the natural map
$$
\Tor_1^{\Z\La P\Ra}(N,P/\fm^{n+1}_P)\to \Tor_1^{\Z\La
P\Ra}(N,P/\fm^n_P)
$$
is injective for every $n\in\N$. Then, a simple induction shows
that, under assumption (c), all these modules vanish. The latter
means that the natural map :
$$
\fm^n_P\otimes_PN\to\fm^n_PN
$$
is an isomorphism, for every $n\in\N$. We consider the commutative
ladder with exact rows :
$$
\xymatrix{ \fm^{n+1}_P\otimes_PN \ar[r] \ar[d] & \fm^n_P\otimes_PN
\ar[r] \ar[d] & (\fm^n_P/\fm^{n+1}_P)\otimes_PN \ar[r] \ar[d] & 0 \\
\fm^{n+1}_PN \ar[r] & \fm^n_PN \ar[r] & \fm^n_PN/\fm^{n+1}_PN \ar[r]
& 0.}
$$
By the foregoing, the two left-most vertical arrows are
isomorphisms, hence the same holds for the right-most, whence (c).

(d)$\Rightarrow$(c): We have to show that the natural map
$u:\fm_P\otimes_PN\to\fm_PN$ is an isomorphism. To this aim, we
consider the $\fm_P$-adic filtrations on these two modules; for the
associated graded modules one gets :
$$
\gr_n(\fm_PN)=\fm^n_PN/\fm^{n+1}_PN \qquad
\gr_n(\fm_P\otimes_PN)=(\fm^n_P/\fm^{n+1}_P)\otimes_PN
$$
for every $n\in\N$; whence maps of $A$-modules :
$$
\gr_n(\fm_P\otimes_PN) \xrightarrow{ \gr_nu }\gr_n(\fm_PN).
$$
which are isomorphism by assumption. To conclude, it suffices to
show :
\begin{claim}\label{cl_loc-flatn-crit}
For every ideal $I\subset P$, the $\fm_P A$-adic filtration is
separated on the $A$-module $I\otimes_PN$.
\end{claim}
\begin{pfclaim} Indeed, notice that the ideal
$I/P^\times\subset P^\sharp$ is finitely generated (proposition
\ref{prop_ideals-in-fg-mon}(ii)), hence the same holds for $I$,
so $I\otimes_PN$ is a finitely generated $A$-module. Then the
contention follows from \cite[Th.8.10]{Mat}.
\end{pfclaim}

(c)$\Rightarrow$(b): We argue by induction on $i$. For $i=1$,
suppose first that $M=P/I$ for some ideal $I\subset P$ (notice that
any such quotient is an integral pointed $P$-module); in this case,
the assertion to prove is that the natural map $v:I\otimes_PN\to IN$
is an isomorphism. However, consider the $\fm_P$-adic filtration on
$P/I$; for the associated graded module we have :
$$
\gr_n(P/I)=(\fm^n_P\cup I)/(\fm^{n+1}\cup I) \qquad \text{for every
$n\in\N$}
$$
and it is easily seen that this is a free pointed $P^\times$-module,
for every $n\in\N$. Hence, by inspecting the long exact
$\Tor$-sequences associated to the exact sequences
$$
0\to\gr_n(P/I)\to P/(\fm^{n+1}_P\cup I)\to P/(\fm^n_P\cup I)\to 0
$$
our assumption (c), together with a simple induction yields :
\set\begin{equation}\label{eq_approx-Tor}
\Tor^{\Z\La P\Ra}_1(N,\Z\La P/(\fm^n_P\cup I)\Ra)
\qquad
\text{for every $n\in\N$.}
\end{equation}
Now, fix $n\in\N$; in light of lemma \ref{lem_log-Artin-Rees}, there
exists $k\geq n$ such that $\fm^k_P\cap I\subset\fm^n_PI$. We deduce
surjective maps of $A$-modules :
$$
I\otimes_PN\to\frac{I}{\fm^k_P\cap I}\otimes_PN\to
\frac{I}{\fm^n_PI}\otimes_PN\isom\frac{I\otimes_PN}{\fm^n_P(I\otimes_PN)}.
$$
On the other hand, \eqref{eq_approx-Tor} says that the natural map
$(\fm^n_P\cup I)\otimes_PN\to(\fm^n_P\cup I)N$ is an isomorphism, so
the same holds for the induced composed map :
$$
\frac{I}{\fm^k_P\cap I}\otimes_PN\isom\frac{\fm^k_P\cup
I}{\fm^k_P}\otimes_PN\to\frac{(\fm^k_P\cup I)N}{\fm^k_PN}.
$$
Consequently, the kernel of $v$ is contained in
$\fm^n_P(I\otimes_PN)$; since $n$ is arbitrary, we are reduced to
showing that the $\fm_P$-adic filtration is separated on
$I\otimes_PN$, which is claim \ref{cl_loc-flatn-crit}.

Next, again for $i=1$, let $M$ be an arbitrary integral $P$-module.
In view of remark \ref{rem_integr-modules}(i), we may assume that
$M$ is finitely generated; moreover, remark
\ref{rem_integr-modules}(ii), together with an easy induction
further reduces to the case where $M$ is cyclic, in which case,
according to remark \ref{rem_integr-modules}(iii), $M$ is of the
form $P/I$ for some ideal $I$, so the proof is complete in this
case.

Finally, suppose $i>1$ and assume that the assertion is already
known for $i-1$. We may similarly reduce to the case where $M=P/I$
for some ideal $I$ as in the foregoing; to conclude, we observe that
:
$$
\Tor_i^{\Z\La P\Ra}(N,\Z\La P/I\Ra)\simeq\Tor_{i-1}^{\Z\La
P\Ra}(N,\Z\La I\Ra).
$$
Since obviously $I$ is an integral $P$-module, the contention
follows.
\end{proof}

As a corollary, we have the following combinatorial going-down
theorem, which is proved in the same way as its commutative algebra
counterpart.

\begin{corollary}\label{cor_going-down-combin}
In the situation of \eqref{subsec_log-local-flat-crit}, assume
that $P^\sharp$ is fine, and that $A$ is $\alpha$-flat. Let
$\fp\subset\fq$ be two prime ideals of $P$, and $\fq'\subset A$
a prime ideal such that $\fq=\alpha^{-1}\fq'$. Then there exists
a prime ideal $\fp'\subset\fq'$ such that $\fp=\alpha^{-1}\fp'$.
\end{corollary}
\begin{proof} Let $\beta:P_\fq\to A_{\fq'}$ be the morphism induced
by $\alpha$; it is easily seen that $A_{\fq'}$ is $\beta$-flat, and
moreover $(P_\fp)^\sharp=(P^\sharp)^\sharp_\fp$ is still fine (lemma
\ref{lem_face}(iv)). Hence we may replace $\alpha$ by $\beta$, which
(in view of \eqref{subsec_sepc-of-monoid}) allows to assume that
$\fq$ (resp. $\fq'$) is the maximal ideal of $P$ (resp. of $A$).
Next, let $P_0:=P/\fp$, $A_0:=A/\fp A$ and denote by
$\alpha_0:P_0\to A_0$ the morphism induced from $\alpha$; it is
easily seen that $A_0$ is $\alpha_0$-flat : for instance, the
natural map
$(\fm^n_P/\fm^{n+1}_P)\otimes_PA_0\to\fm^n_PA_0/\fm_P^{n+1}A_0$ is
of the type $f\otimes_AA_0$, where $f$ is the map in proposition
\ref{prop_flat-vrit}(c), thus if the latter is bijective, so is the
former. Moreover $P_0^\sharp$ is a quotient of $P^\sharp$,
hence it is again fine. Therefore we may replace $P$ by $P_0$ and
$A$ by $A_0$, which allows to further assume that $\fp=\{0\}$, and
it suffices to show that there exists a prime ideal $\fq'\subset A$,
such that $\alpha^{-1}\fq'=\{0\}$. Set $\Sigma:=P\setminus\{0\}$; it
is easily seen that the natural morphism $P\to\Sigma^{-1}P$ is
injective; moreover, its cokernel $C$ (in the category of pointed
$P$-modules) is integral, so that $\Tor_1^{\Z\La P\Ra}(A,\Z\La
C\Ra)=0$ by assumption. It follows that the localization map
$A\to\Sigma^{-1}A$ is injective, especially $\Sigma^{-1}A\neq\{0\}$,
and therefore it contains a prime ideal $\fq''$. The prime ideal
$\fq':=\fq''\cap A$ will do.
\end{proof}

\begin{lemma}\label{lem_invariance-of-flatness}
Let $A$ be a noetherian local ring, $N$ an $A$-module of finite
type, $\alpha:P\to A$ and $\beta:Q\to A$ two morphisms of pointed
monoids, with $P^\sharp$ and $Q^\sharp$ both fine. Suppose that
$\alpha$ and $\beta$ induce the same constant log structure on
$\Spec\,A$ (see \eqref{subsec_Konstant}). Then $N$ is $\alpha$-flat
if and only if it is $\beta$-flat.
\end{lemma}
\begin{proof} Let $\xi$ be a $\tau$-point localized at the closed
point of $X:=\Spec\,A$, set $B:=\cO_{\!X,\xi}$ and let $\phi:A\to B$
be the natural map. Let also $M$ be the push-out of the diagram of
monoids $P\leftarrow(\phi\circ\alpha)^{-1}B^\times\to B^\times$
deduced from $\alpha$; then $M\simeq P_{\!X,\xi}^{\log}$, the stalk
at the point $\xi$ of the constant log structure on $X_\tau$
associated to $\alpha$. Since $\phi$ is faithfully flat, it is
easily seen that $N$ is $\alpha$-flat if and only if $N\otimes_AB$
is $\phi\circ\alpha$-flat.

Hence we may replace $A$ by $B$, $\alpha$ by $\phi\circ\alpha$, $N$
by $N\otimes_AB$, and $Q$ by $M$, after which we may assume that
$Q=P\amalg_{\alpha^{-1}(A^\times)}A^\times$; especially, there
exists a morphism of monoids $\gamma:P\to Q$ such that
$\alpha=\beta\circ\gamma$, and moreover $\beta$ is a local morphism.

Next, set $S:=\alpha^{-1}(A^\times)$; clearly $\gamma$ extends to a
morphism of monoids $\gamma':S^{-1}P\to Q$, and $\alpha$ and
$\beta\circ\gamma'$ induce the same constant log structure on
$X_\tau$. Arguing as in the proof of proposition
\ref{prop_flat-vrit}, we see that $N$ is $P$-flat if and only if it
is $S^{-1}P$-flat. Hence, we may replace $P$ by $S^{-1}P$, which
allows to assume that $\gamma$ induces an isomorphism
$P\amalg_{P^\times}A^\times\isom Q$, therefore also an isomorphism
$P^\sharp\isom Q^\sharp$. The latter implies that
$\fm_Q=\fm_PQ$; moreover, notice that the morphism of monoids
$P^\times\to A^\times$ is faithfully flat, so the natural map :
$$
\Tor_1^{\Z\La P\Ra}(N,\Z\La P/\fm_P\Ra)\to\Tor_1^{\Z\La
Q\Ra}(N,\Z\La (P/\fm_P)\otimes_PQ\Ra)\to\Tor^{\Z\La Q\Ra}_1(N,\Z\La
Q/\fm_Q\Ra)
$$
is an isomorphism. The assertion follows.
\end{proof}

\begin{lemma}\label{lem_subsheaf-if-flat}
Let $M$ be an integral monoid, $A$ a ring, $\phi:M\to A$ a morphism
of monoids, and set $S:=\Spec\,A$. Suppose that $A$ is $\phi$-flat.
Then the log structure $(M,\phi)^{\log}_S$ on $S_\tau$ is the subsheaf
of monoids of $\cO_{\!S}$ generated by $\cO^\times_{\!S}$ and the
image of $M$.
\end{lemma}
\begin{proof} To ease notation, set $\underline
M:=(M,\phi)^{\log}_S$; let $\xi$ be any $\tau$-point of $S$. Then
the stalk $\underline M{}_\xi$ is the push-out of the diagram :
$\cO^\times_{\!S,\xi}\leftarrow
\phi_\xi^{-1}(\cO_{\!S,\xi}^\times)\to M$ where
$\phi_\xi:M\to\cO_{\!S,\xi}$ is deduced from $\phi$. Hence
$\underline M{}_\xi$ is generated by $\cO_{\!S,\xi}^\times$ and the
image of $M$, and it remains only to show that the structure map
$\underline M{}_\xi\to\cO_{\!S,\xi}$ is injective. Therefore, let
$a,b\in M$ and $u,v\in\cO_{\!S,\xi}^\times$ such that :
\set\begin{equation}\label{eq_lead-to-equo}
\phi_\xi(a)\cdot u=\phi_\xi(b)\cdot v.
\end{equation}
We come down to showing :
\begin{claim} There exist $c,d\in M$ such that :
$$
\phi_\xi(c),\phi_\xi(d)\in\cO_{\!S,\xi}^\times
\qquad
ac=bd
\qquad
\phi_\xi(c)\cdot v=\phi_\xi(d)\cdot u.
$$
\end{claim}
\begin{pfclaim}[] Let $\fm_\xi\subset\cO_{\!S,\xi}$ be the maximal
ideal, and set $\fp:=\phi_\xi^{-1}(\fm_\xi)$, so that $\phi_\xi$
extends to a local morphism $\phi_\fp:M_\fp\to\cO_{\!S,\xi}$. Since
$\cO_{\!S,\xi}$ is a flat $A$-algebra, $\cO_{\!S,\xi}$ is
$\phi$-flat, and consequently it is faithfully $\phi_\fp$-flat
(lemma \ref{lem_faithful-phi-flat}). Then, assumption
\eqref{eq_lead-to-equo} leads to the identity:
$$
aM_\fp\otimes_{M_\fp}\cO_{\!S,\xi}=\phi_\xi(a)\cdot\cO_{\!S,\xi}=
\phi_\xi(b)\cdot\cO_{\!S,\xi}=bM_\fp\otimes_{M_\fp}\cO_{\!S,\xi}
$$
whence $aM_\fp=bM_\fp$, by faithful $\phi_\fp$-flatness. It follows
that there exist $x,y\in M_\fp$ such that $ax=b$ and $by=a$, hence
$axy=a$, which implies that $xy=1$, since $M_\fp$ is an integral
module. The latter means that there exist $c,d\in M\!\setminus\!\fp$
such that $ac=bd$ in $M$. We deduce easily that
$$
\phi_\xi(a)\cdot\phi_\xi(c)\cdot v=\phi_\xi(a)\cdot\phi_\xi(d)\cdot
u.
$$
Thus, to complete the proof, it suffices to show that $\phi_\xi(a)$
is regular in $\cO_{\!S,\xi}$. However, the morphism of $M$-modules
$\mu_a:M\to M$ : $m\mapsto am$ (for all $m\in M$) is injective,
hence the same holds for the map
$\mu_a\otimes_M\cO_{\!S,\xi}:\cO_{\!S,\xi}\to\cO_{\!S,\xi}$, which
is just multiplication by $\phi(a)$.
\end{pfclaim}
\end{proof}

\sset\subsubsection{}\label{subsec_define-new-dim}
Let $(X,\underline M)$ be a locally noetherian log scheme, with
coherent log structure (on the site $X_\tau$, see
\eqref{subsec_special-schs}), and let $\xi$ be any $\tau$-point
of $X$. We denote by $I(\xi,\underline M)\subset\cO_{\!X,\xi}$ the
ideal generated by the image of the maximal ideal of
$\underline M{}_\xi$, and we set :
$$
d(\xi,\underline M):=
\dim\cO_{\!X,\xi}/I(\xi,\underline M)+\dim\underline M{}_\xi.
$$

\begin{lemma}\label{lem_we-got-a-situation}
In the situation of \eqref{subsec_define-new-dim}, suppose
furthermore that $(X,\underline M)$ is a fs log scheme.
Then we have the inequality :
$$
\dim\cO_{X,\xi}\leq d(\xi,\underline M).
$$
\end{lemma}
\begin{proof} According to corollary \ref{cor_sharpy}(i), there exist
a neighborhood $U\to X$ of $\xi$ in $X_\tau$, and a fine and
saturated chart $\alpha:P_U\to\underline M_{|U}$, which is sharp at
the point $\xi$. Especially, $P\simeq\underline M{}_\xi^\sharp$,
therefore $\dim P=\dim\underline M{}_\xi$ (corollary
\ref{cor_consequent}(ii)). Notice that $\cO_{\!X,\xi}$ is a
noetherian local ring (this is obvious for $\tau=\text{Zar}$, and
follows from \cite[Ch.IV, Prop.18.8.8(iv)]{EGA4} for
$\tau=\text{\'et}$), hence to conclude it suffices to apply
corollary \ref{cor_new-dim}(i) to the induced map of monoids
$P\to\cO_{\!X,\xi}$.
\end{proof}

\begin{definition}\label{def_log-regular}
Let $(X,\underline M)$ be a locally noetherian fs log scheme, $\xi$
a $\tau$-point of $X$.
\begin{enumerate}
\item
We say that $(X,\underline M)$ is {\em regular at the point\/} $\xi$,
if the following holds :
\begin{enumerate}
\item
the inequality of lemma \eqref{lem_we-got-a-situation} is actually
an equality, and
\item
the local ring $\cO_{\!X,\xi}/I(\xi,\underline M)$ is regular.
\end{enumerate}
\item
We denote by $(X,\underline M)_\reg\subset X$ the set of points $x$
such that $(X,\underline M)$ is regular at any (hence all)
$\tau$-points of $X$ localized at $x$.
\item
We say that $(X,\underline M)$ is {\em regular}, if
$(X,\underline M)_\reg=X$.
\item
Suppose that $K$ is a field, and $X$ a $K$-scheme. We say that
the $K$-log scheme $(X,\underline M)$ is {\em geometrically regular},
if $E\times_K(X,\underline M)$ is regular, for every field
extension $E$ of $K$.
\end{enumerate}
\end{definition}

\begin{remark}\label{rem_pointed-reg}
(i)\ \
Certain constructions produce log structures
$\underline M\to\cO_{\!X}$ that are morphisms of pointed
monoids. It is then useful to extend the notion of regularity
to such log structures. We shall say that $(X,\underline M)$
is a {\em pointed regular\/} log scheme, if there exists a
log structure $\underline N$ on $X$, such that
$\underline M=\underline N{}_\circ$ (notation of
\eqref{subsec_pointed-log-topos}), and $(X,\underline N)$
is a regular log scheme.

(ii)\ \
Likewise, if $K$ is a field, $X$ a $K$-scheme, and
$\underline M$ a log structure on $X$, we shall say that
the $K$-log scheme $(X,\underline M)$ is {\em geometrically
pointed regular}, if $\underline M=\underline N{}_\circ$
for some log structure $\underline N$ on $X$, such that
$(X,\underline N)$ is geometrically regular.
\end{remark}

\sset\subsubsection{}\label{subsec_fix-some-nota}
Let $(X,\underline M)$ be a locally noetherian fs log scheme,
$\xi$ a $\tau$-point of $X$, and $\cO^\wedge_{\!X,\xi}$ the
completion of $\cO_{\!X,\xi}$. The next result is the logarithmic
version of the classical characterization of complete regular
local rings (\cite[Th.29.7 and Th.29.8]{Mat}).

\begin{theorem}\label{th_charact-log-regular}
With the notation of \eqref{subsec_fix-some-nota}, the log scheme
$(X,\underline M)$ is regular at the point $\xi$ if and only if
there exist :
\begin{enumerate}
\alphaenu
\item
a complete regular local ring $(R,\fm_R)$, and a local ring homomorphism
$R\to\cO^\wedge_{\!X,\xi}$;
\item
a fine and saturated chart $P_{X(\xi)}\to\underline M(\xi)$ which is
sharp at the closed point $\xi$ of $X(\xi)$, such that the induced
continuous ring homomorphism
$$
R[[P]]\to\cO^\wedge_{\!X,\xi}
$$
is an isomorphism if $\cO_{\!X,\xi}$ contains a field, and otherwise
it is a surjection, with kernel generated by an element $\theta\in
R[[P]]$ whose constant term lies in $\fm_R\!\setminus\!\fm_R^2$.
\end{enumerate}
\end{theorem}
\begin{proof} Suppose first that (a) and (b) hold.  If $R$ contains
a field, then it follows that
$$
\dim\cO_{\!X,\xi}=\dim\cO_{\!X,\xi}^\wedge=\dim R+\dim P
$$
by corollary \ref{cor_new-dim}(ii); moreover, in this case
$I(\xi,\underline M)=\fm_P\cO_{\!X,\xi}$, hence
$\cO_{\!X,\xi}/I(\xi,\underline
M)=\cO_{\!X,\xi}/\fm_P\cO_{\!X,\xi}$, whose completion is
$\cO^\wedge_{\!X,\xi}/\fm_P\cO^\wedge_{\!X,\xi}\simeq R$ so that
$\cO_{\!X,\xi}/I(\xi,\underline M)$ is regular (\cite[Ch.0,
Prop.17.3.3(i)]{EGAIV}). Furthermore,
$$
\dim R=\dim\cO^\wedge_{\!X,\xi}/\fm_P\cO^\wedge_{\!X,\xi}
=\dim\cO_{\!X,\xi}/I(\xi,\underline M)
$$
(\cite[Th.15.1]{Mat}). Hence $(X,\underline M)$ is regular at the
point $\xi$. If $R$ does not contain a field, we obtain
$\dim\cO_{\!X,\xi}=\dim R+\dim P-1$. On the other hand, let
$\theta_0$ be the image of $\theta$ in $\fm_R$; then we have
$\cO^\wedge_{\!X,\xi}/\fm_P\cO^\wedge_{\!X,\xi}\simeq R/\theta_0R$,
which is regular of dimension $\dim\,R-1$, and again we invoke
\cite[Ch.0, Prop.17.3.3(i)]{EGAIV} to see that $(X,\underline M)$ is
regular at $\xi$.

Conversely, suppose that $(X,\underline M)$ is regular at $\xi$.
Suppose first that $\cO_{\!X,\xi}$ contains a field; then we may
find a field $k\subset\cO_{\!X,\xi}^\wedge$ mapping isomorphically
to the residue field of $\cO_{\!X,\xi}^\wedge$
(\cite[Th.28.3]{Mat}). Pick a sequence $(t_1,\dots,t_r)$ of
elements of $\cO_{\!X,\xi}$ whose image in the regular local ring
$\cO_{\!X,\xi}/I(\xi,\underline M)$ forms a regular system of
parameters. Let also $P$ a fine saturated monoid for which there
exist a neighborhood $U\to X$ of $\xi$ and a chart
$P_U\to\underline M_{|U}$, sharp at the point $\xi$. There follows
a map of monoids $\alpha:P\to\cO_{\!X,\xi}$, and necessarily the
image of $\fm_P$ lies in the maximal ideal of $\cO_{\!X,\xi}$, and
generates $I(\xi,\underline M)$. We deduce a continuous ring
homomorphism
$$
k[[P\times\N^{\oplus r}]]\to\cO^\wedge_{\!X,\xi}
$$
which extends $\alpha$, and which maps the generators
$T_1,\dots,T_r$ of $\N^{\oplus r}$ onto respectively
$t_1,\dots,t_r$. This map is clearly surjective, and by comparing
dimensions (using corollary \ref{cor_new-dim}(ii)) one sees that it
is an isomorphism. Then the theorem holds in this case, with
$R:=k[[\N^{\oplus r}]]$.

Next, if $\cO_{\!X,\xi}$ does not contain a field, then its residue
characteristic is a positive integer $p$, and we may find a complete
discrete valuation ring $V\subset\cO^\wedge_{\!X,\xi}$ whose maximal
ideal is $pV$, and such that $V/pV$ maps isomorphically onto the
residue field of $\cO^\wedge_{\!X,\xi}$ (\cite[Th.29.3]{Mat}).
Again, we choose a morphism of monoids $\alpha:P\to\cO_{\!X,\xi}$ as
in the foregoing, and a sequence $(t_1,\dots,t_r)$ of elements of
$\cO_{\!X,\xi}$ lifting a regular system of parameters for
$\cO_{\!X,\xi}/I(\xi,\underline M)$, by means of which we define a
continuous ring homomorphism
$$
\phi:V[[P\times\N^{\oplus r}]]\to\cO^\wedge_{\!X,\xi}
$$
as in the previous case. Again, it is clear that $\phi$ is
surjective. The image of the ideal $J$ generated by the maximal
ideal of $P\times\N^{\oplus r}$ is the maximal ideal of
$\cO^\wedge_{\!X,\xi}$; in particular, there exists $x\in J$ such
that $\theta:=p-x$ lies in $\Ker\,\phi$. If we let $R:=V[[\N^{\oplus
r}]]$, it is clear that $\theta\in\fm_R\setminus\fm_R^2$.

\begin{claim}\label{cl_separ}
Let $A$ be a ring, $\pi$ a regular element of $A$ such that $A/\pi
A$ is an integral domain, $P$ a fine and sharp monoid. Let also
$\theta$ be an element of $A[[P]]$ whose constant term is $\pi$
(see \eqref{subsec_constant-term}). Then $A[[P]]/(\theta)$ is an
integral domain.
\end{claim}
\begin{pfclaim} To ease notation, set $A_0:=A/\pi A$, $B:=A[[P]]$
and $C:=B/\theta B$. Choose a decomposition \eqref{eq_decomponi},
and set $\Fil^\gamma_nB:=\prod_{i\geq n}A[[\gr^\gamma_iP]]$ for
every $n\in\N$. $\Fil^\gamma_\bullet B$ is a separated filtration by
ideals of $B$, and we may consider the induced filtration
$\Fil^\gamma_\bullet C$ on $C$. First, we remark that
$\Fil^\gamma_\bullet C$ is also separated. This comes down to
checking that
$$
\bigcap_{n\geq 0}\theta B+\Fil^\gamma_nB=\theta B.
$$
To this aim, suppose that, for a given $x\in B$ we have identities
of the type $x=\theta y_n+z_n$, with $y_n\in B$ and
$z_n\in\Fil^\gamma_n B$ for every $n\in\N$. Then, since $\pi$ is
regular, an easy induction shows that
$\gr^\gamma_i(y_n)=\gr^\gamma_i(y_m)$ whenever $n,m>i$, and moreover
$x=\theta\cdot\sum_{i\in\N}\gr^\gamma_i(y_{i+1})$, which shows the
contention. It follows that, in order to show that $C$ is a domain,
it suffices to show that the same holds for the graded ring
$\gr^\gamma_\bullet C$ associated to $\Fil^\gamma_\bullet C$.
However, notice that :
$$
\Fil^\gamma_n B\cap\theta B=\theta\cdot\Fil^\gamma_n B \qquad
\text{for every $n\in\N$}.
$$
(Indeed, this follows easily from the fact that $\pi$ is a regular
element : the verification shall be left to the reader). Hence, we
may compute :
$$
\gr^\gamma_n C=\frac{\Fil^\gamma_nB+\theta
B}{\Fil^\gamma_{n+1}+\theta B}\simeq\frac{\Fil^\gamma_n
B}{\Fil^\gamma_{n+1}+\theta\Fil^\gamma_n}\simeq A[\gr^\gamma_n
P]/\theta A[\gr^\gamma_n P]\simeq A_0[\gr^\gamma_n P].
$$
Thus, $\gr^\gamma_\bullet A\simeq A_0[P]$, which is a domain, since
by assumption $A_0$ is a domain.
\end{pfclaim}

From claim \ref{cl_separ}(ii) we deduce that $R[[P]]/(\theta)$ is an
integral domain. Then, again by comparing dimensions, we see that
$\phi$ factors through an isomorphism
$R[[P]]/(\theta)\isom\cO^\wedge_{\!X,\xi}$.
\end{proof}

\begin{remark}\label{rem_more-precisely}
Resume the notation of \eqref{subsec_fix-some-nota}, and
suppose that $\cO_{\!X,\xi}/I(\xi,\underline M)$ is a regular local
ring. Let $P_{X(\xi)}\to\underline M(\xi)$ be a chart as in theorem
\ref{th_charact-log-regular}(b), and $\alpha:P\to\cO_{\!X,\xi}$ the
corresponding morphism of monoids. Moreover, if $\cO_{\!X,\xi}$
contains a field, let $V$ denote a coefficient field of
$\cO^\wedge_{\!X,\xi}$, and otherwise, let $V$ be a complete
discrete valuation ring whose maximal ideal is generated by $p$, the
residue characteristic of $\cO_{\!X,\xi}$. In either case, pick a
ring homomorphism $V\to\cO^\wedge_{\!X,\xi}$ inducing an isomorphism
of $V/pV$ onto the residue field of $\cO_{\!X,\xi}$. Let as well
$t_1,\dots,t_r\in\cO_{\!X,\xi}$ be any sequence of elements whose
image in $\cO_{\!X,\xi}/I(\xi,\underline M)$ forms a regular system
of parameters, and extend the map $\alpha$ to a morphism of monoids
$P\times\N^{\oplus r}\to\cO_{\!X,\xi}$, by the rule : $e_i\mapsto
t_i$, where $e_1,\dots,e_r$ is the natural basis of $\N^{\oplus r}$.
Then by inspecting the proof of theorem
\ref{th_charact-log-regular}, we see that the induced continuous
ring homomorphism :
$$
V[[P\times\N^{\oplus r}]]\to\cO^\wedge_{\!X,\xi}
$$
is always surjective, and if $V$ is not a field, its kernel contains
an element $\theta$ whose constant term in $V$ is $\theta_0=p$. If
$V$ is a field (resp. if $V$ is a discrete valuation ring) then
$(X,\underline M)$ is regular at the point $\xi$, if and only if
this map is an isomorphism (resp. if and only if the kernel of this
map is generated by $\theta$).
\end{remark}

\begin{corollary}\label{cor_normal-and-CM}
Let $(X,\underline M)$ be a regular log scheme. Then the scheme $X$
is normal and Cohen-Macaulay.
\end{corollary}
\begin{proof} Let $x\in X$ be any point, and $\xi$ a $\tau$-point
localized at $x$; we have to show that $\cO_{\!X,x}$ is
Cohen-Macaulay; in light of \cite[Ch.IV, Cor.18.8.13(a)]{EGA4} (when
$\tau=\et$), it suffices to show that the same holds for
$\cO_{\!X,\xi}$. Then \cite[Th.17.5]{Mat} further reduces to showing
that the completion $\cO^\wedge_{\!X,\xi}$ is Cohen-Macaulay; the
latter follows easily from theorems \ref{th_charact-log-regular} and
\ref{th_Hochster}(i). Next, in order to prove that $X$ is normal, it
suffices to show that $\cO_{\!X,x}$ is regular, whenever $x$ has
codimension one in $X$ (\cite[Th.23.8]{Mat}). Again, by \cite[Ch.IV,
Cor.18.8.13(c)]{EGA4} (when $\tau=\et$) and \cite[Ch.0,
Prop.17.1.5]{EGAIV}, we reduce to showing that
$\cO_{\!X,\xi}^\wedge$ is regular for such a point $x$. However, for
a point of codimension one we have $r:=\dim\underline M_\xi\leq 1$.
If $r=0$, then $I(\xi,\underline M)=\{0\}$, hence $\cO_{\!X,\xi}$ is
regular. Lastly, if $r=1$, we see that $\underline M_\xi^\sharp\simeq\N$
(theorem \ref{th_structure-of-satu}(iii));
consequently, there exist a regular local ring $R$ and an
isomorphism $\cO^\wedge_{\!X,\xi}\simeq R[[\N]]/(\theta)$, where
$\theta=0$ if $\cO_{\!X,\xi}$ contains a field, and otherwise the
constant term of $\theta$ lies in $\fm_R\setminus\fm^2_R$. Since
$R[[\N]]$ is again regular, the assertion follows in either case.
\end{proof}

\begin{corollary}\label{cor_regular-immers}
Let $i:(X',\underline M')\to(X,\underline M)$ be an exact closed
immersion of regular log schemes (see definition
{\em\ref{def_Exact-immersion}(i)}). Then the underlying morphism of
schemes $X'\to X$ is a regular closed immersion.
\end{corollary}
\begin{proof} Let $\xi$ be a $\tau$-point of $X$, and denote by
$J$ the kernel of $i^\natural_\xi:\cO_{\!X,\xi}\to\cO_{\!X',\xi}$.
To ease notation, let as well $A:=\cO_{\!X,\xi}/I(\xi,\underline M)$
and $A':=\cO_{\!X',\xi}/I(\xi,\underline M')$ Since $\log
i:i^*\underline M\to\underline M'$ is an isomorphism,
$i^\natural_\xi$ induces an isomorphism :
$$
\cO_{\!X,\xi}/(J+I(\xi,\underline M))\isom A'.
$$
Since $A$ and $A'$ are regular, there exists a sequence of elements
$t_1,\dots,t_k\in J$, whose image in $A$ forms the beginning of a
regular system of parameters and generate the kernel of the induced
map $A\to A'$ (\cite[Ch.0, Cor.17.1.9]{EGAIV}). Extend this sequence
by suitable elements of $\cO_{\!X,\xi}$, to obtain a sequence
$(t_1,\dots,t_r)$ whose image in $A$ is a regular system of
parameters. We deduce a surjection $\phi:V[[P\times\N^{\oplus
r}]]\to\cO^\wedge_{\!X,\xi}$ as in remark \ref{rem_more-precisely},
where $V$ is either a field or a complete discrete valuation ring.
Denote by $(e_1,\dots,e_r)$ the natural basis of $\N^{\oplus r}$.
Now, suppose first that $V$ is a complete discrete valuation ring;
then $\Ker\,\phi$ is generated by an element $\theta\in
V[[P\times\N^{\oplus r}]]$, whose constant term is a uniformizer in
$V$; we deduce easily from claim \ref{cl_separ} that
$(e_1,\dots,e_r,\theta)$ is a regular sequence in
$V[[P\times\N^{\oplus r}]]$; hence the same holds for the sequence
$(\theta,e_1,\dots,e_r)$, in view of \cite[p.127, Cor.]{Mat}. This
implies that $(t_1,\dots,t_r)$ is a regular sequence in
$\cO^\wedge_{\!X,\xi}$, hence $(t_1,\dots,t_k)$ is a regular
sequence in $\cO_{\!X,\xi}$, which is the contention. The case where
$V$ is a field is analogous, though simpler : the details shall be
left to the reader.
\end{proof}

\begin{remark}\label{rem_simplex-is-normal-cross}
In the situation of remark \ref{rem_more-precisely}, suppose that
$P$ is a free monoid of finite rank $d$, let
$\{f_1,\dots,f_d\}\subset\cO_{\!X,\xi}$ be the image of the (unique)
basis of $P$, and for every $i\leq d$ let $Z_i\subset X$ be the
zero locus of $f_i$; then $\bigcup_{i=1}^d Z_i$ is a strict normal
crossings divisor in the sense of example \ref{ex_norm-cross}. The
proof is analogous to that of corollary \ref{cor_regular-immers} :
the sequence $(f_1,\dots,f_d,\theta)$ is regular in
$V[[P\oplus\N^{\oplus r}]]$, hence also its permutation
$(\theta,f_1,\dots,f_d)$ is regular, whence the claim (the details
shall be left to the reader).
\end{remark}

\begin{proposition}\label{prop_second-crit}
Resume the situation of \eqref{subsec_fix-some-nota}. We have :
\begin{enumerate}
\item
The log scheme $(X,\underline M)$ is regular at the $\tau$-point
$\xi$ if and only if the following two conditions hold :
\begin{enumerate}
\item
The ring $\cO_{\!X,\xi}/I(\xi,\underline M)$ is regular.
\item
There exists a morphism of monoids $P\to\cO_{\!X,\xi}$ from a fine
monoid $P$, whose associated constant log structure on $X(\xi)$ is
the same as $\underline M(\xi)$, and such that $\cO_{\!X,\xi}$ is
$P$-flat.
\end{enumerate}
\item
Moreover, if the equivalent conditions of {\em (i)} hold, and
$Q_X\to\underline M(\xi)$ is any fine chart, then $\cO_{\!X,\xi}$
is $Q$-flat, for the induced map of monoids $Q\to\cO_{\!X,\xi}$. 
\end{enumerate}
\end{proposition}
\begin{proof}(i): Suppose first that $(X,\underline M)$ is regular at
the point $\xi$; then by definition (a) holds. Next, we may find a
fine monoid $P$ and a local morphism $\alpha:P\to A:=\cO_{\!X,\xi}$
whose associated constant log structure on $X(\xi)$ is the same as
$\underline M(\xi)$, and a regular local ring $R$, with a ring
homomorphism $R\to A$ such that the induced continuous map
$R[[P]]\to A^\wedge$ fulfills condition (b) of theorem
\ref{th_charact-log-regular} (where $A^\wedge$ is the completion of
$A$). Since the completion map $\phi:A\to A^\wedge$ is faithfully
flat, $A$ is $\alpha$-flat if and only if $A^\wedge$ is
$\phi\circ\alpha$-flat; in light of proposition \ref{rem_whenever-phi}(i),
it then suffices to show that, for every ideal $I\subset P$, the
natural map
$$
A^\wedge\derotimes_PP/I\to A^\wedge\otimes_PP/I
$$
is an isomorphism in $\sD^-(A^\wedge\Mod)$ (notation of
\eqref{subsec_restrict-scalar-mon}). To this aim we remark :
\begin{claim}\label{cl_WinEDT}
For any ideal $I\subset P$, the natural morphism
$$
R[P]\derotimes_PP/I\to R\La P/I\Ra
$$
is an isomorphism in $\sD^-(R[P]\Mod)$.
\end{claim}
\begin{pfclaim} We consider the change of ring spectral
sequence for $\Tor$ :
$$
E^2_{ij}:=\Tor^{\Z[P]}_i(\Tor^\Z_j(R,\Z[P]),\Z\La P/I\Ra)
\Rightarrow\Tor_{i+j}^\Z(R,\Z\La P/I\Ra)
$$
(\cite[Th.5.6.6]{We}). Clearly $E^2_{ij}=0$ whenever $j>0$, whence
isomorphisms :
$$
\Tor^{\Z\La P\Ra}_i(R\La P\Ra,\Z\La P/I\Ra)\isom\Tor_i^\Z(R,\Z\La
P/I\Ra)
$$
for every $i\in\N$. The claim follows easily.
\end{pfclaim}

In light of claim \ref{cl_WinEDT} we deduce natural isomorphisms :
$$
A^\wedge\derotimes_PP/I\isom
A^\wedge\derotimes_{R[P]}(R[P]\derotimes_PP/I)\isom
A^\wedge\derotimes_{R[P]}R\La P/I\Ra
$$
in $\sD^-(A^\wedge\Mod)$. Now, if $A$ contains a field, we have
$A^\wedge\simeq R[[P]]$, which is a flat $R[P]$-algebra
(\cite[Th.8.8]{Mat}), and the contention follows. If $A$ does not
contain a field, the complex
$$
0\to R[[P]]\xrightarrow{ \theta }R[[P]]\to A^\wedge\to 0
$$
is a $R[P]$-flat resolution of $A^\wedge$. Since
$R[[P]]/IR[[P]]\simeq\lim_{n\in\N}\,R\La P/(I\cup\fm^n_P)\Ra$, we
come down to the following :
\begin{claim} Let $\theta\in R[P]$ be any element whose constant term
$\theta_0$ is a regular element in $R$. Then, for every ideal
$I\subset P$, the image of $\theta$ in $R\La P/I\Ra$ is a regular
element.
\end{claim}
\begin{pfclaim} In view of \eqref{subsec_noncanonick}, $P$ can be
regarded as a graded monoid $P=\coprod_{n\in\N}P_n$, so $R[P]$ is a
graded algebra, and $I=\coprod_{n\in\N}(I\cap P_n)$ is a graded
ideal. Thus, $\La P/I\Ra$ is a graded $R$-algebra as well, and the
claim follows easily (details left to the reader).
\end{pfclaim}

Conversely, suppose that conditions (a) and (b) hold. By virtue
of lemma \ref{lem_invariance-of-flatness}, we may assume that
$P$ is fine, sharp and saturated, and that the map $P\to\cO_{\!X,\xi}$
is local. According to remark \ref{rem_more-precisely}, we may
find a regular local ring $R$ of dimension $\leq 1+\dim A/\fm_P A$,
and a surjective ring homomorphism :
$$
\phi:R[[P]]\to A^\wedge
$$
Suppose first that $\dim R=1+\dim A/\fm_PA$; then $\Ker\,\phi$
contains an element $\theta$ whose constant term $\theta_0$ lies in
$\fm_R\setminus\fm_R^2$, and $\phi$ induces an isomorphism
$R_0:=R/\theta_0R\isom A^\wedge/\fm_PA^\wedge$. We have to show that
$\phi$ induces an isomorphism $\bar\phi:R[[P]]/(\theta)\isom
A^\wedge$. To this aim, we consider the $\fm_P$-adic filtrations on
these rings; for the associated graded rings we have :
$$
\gr_nR[[P]]/(\theta)\simeq R_0\La\fm_P^n/\fm_P^{n+1}\Ra \qquad
\gr_nA^\wedge=\fm_P^nA^\wedge/\fm^{n+1}_PA^\wedge.
$$
Hence $\gr_n\bar\phi$ is the natural map
$A^\wedge/\fm_PA^\wedge\otimes_P(\fm^n_P/\fm^{n+1}_P)\to
\fm_P^nA^\wedge/\fm^{n+1}_PA^\wedge$, and the latter is an
isomorphism, since $A^\wedge$ is $P$-flat. The assertion follows in
this case. The remaining case where $\dim R=\dim A/\fm_PA$ is
similar, but easier, so shall be left to the reader, as an exercise.

(ii) follows immediately from (i) and lemma
\ref{lem_invariance-of-flatness}.
\end{proof}

\begin{corollary}\label{cor_more-precisely}
Let $(X,\underline M)$ be a log scheme, $\xi$ a $\tau$-point
of $X$, and suppose that $(X,\underline M)$ is regular at $\xi$.
Then the following conditions are equivalent :
\begin{enumerate}
\alphaenu
\item
$\cO_{\!X,\xi}$ is a regular local ring.
\item
$\underline M{}_\xi^\sharp$ is a free monoid of finite rank.
\end{enumerate}
\end{corollary}
\begin{proof} (b)$\Rightarrow$(a) : Indeed, it suffices to show
that $\cO^\wedge_{\!X,\xi}$ is regular
(\cite[Ch.0, Prop.17.3.3(i)]{EGAIV}); by theorem
\ref{th_charact-log-regular}, the latter is isomorphic to
either $R[[P]]$ or $R[[P]]/\theta R[[P]]$, where $R$ is a
regular local ring, $P:=\underline M{}_\xi$, and the constant
term of $\theta$ lies in $\fm_R\setminus\fm_R^2$. Since, by
assumption, $P$ is a free monoid of finite rank, it is easily
seen that rings of the latter kind are regular (\cite[Ch.0,
Cor.17.1.8]{EGAIV}).

(a)$\Rightarrow$(b) : By proposition \ref{prop_second-crit}(ii),
$R:=\cO_{\!X,\xi}$ is $P$-flat, for a morphism
$P:=\underline M{}_\xi^\sharp\to R$ that induces the log
structure $\underline M(\xi)$ on $X(\xi)$. It follows that
$$
\fm_PR/\fm^2_PR=(\fm_P/\fm_P^2)\otimes_PR=R^{\oplus r}
$$
where $r:=\rk^\circ_{P^\times}\fm_P/\fm_P^2$ is the cardinality
of $\fm_P\setminus\fm_P^2$. In view of \cite[Ch.IV, Prop.16.9.3,
Cor.16.9.4, Cor.19.1.2]{EGA4}, we deduce that the image of
$\fm_P\setminus\fm_P^2$ is a regular sequence of $R$ of length
$r$, hence $\dim P=\dim R-\dim R/\fm_PR=r$, since $(X,\underline M)$
is regular at $\xi$. Then the assertion follows from proposition
\ref{prop_regul-critt}.
\end{proof}

\begin{corollary}\label{cor_same-height}
Let $(X,\underline M)$ be a regular log scheme, $\xi$ any
$\tau$-point of $X$, and $\fp\subset\underline M_\xi$ an
ideal. We have :
\begin{enumerate}
\item
If $\fp$ is a prime ideal, $\fp\cO_{\!X,\xi}$ is a prime ideal
of $\cO_{\!X,\xi}$, and $\hgt\,\fp=\hgt\,\fp\cO_{\!X,\xi}$.
\item
If $\fp$ is a radical ideal, $\fp\cO_{\!X,\xi}$ is a radical ideal
of $\cO_{\!X,\xi}$.
\end{enumerate}
\end{corollary}
\begin{proof} To ease notation, set $A:=\cO_{\!X,\xi}$. Pick a
chart $\beta:P_{X(\xi)}\to\underline M(\xi)$ with $P$ fine,
sharp and saturated (corollary \ref{cor_sharpy}(i)); by
proposition \ref{prop_second-crit}(ii), the ring $A$ is
$P$-flat for the resulting map $P\to A$. Let
$\fq:=\beta^{-1}\fp\subset P$; then $\fp A=\fq A$.

(i): Since $\fp$ is a prime ideal, $\fq$ is a prime ideal of $P$,
and in order to prove that $\fp A$ is a prime ideal, it suffices
therefore to show that the completion $(A/\fq A)^\wedge$ of
$A/\fq A$ is an integral domain. However, remark
\ref{rem_more-precisely} implies that
$(A/\fq A)^\wedge\simeq B/\theta B$, where
$B:=R[[P\!\setminus\!\fq]]$, with $(R,\fm_R)$ a regular local ring,
and $\theta$ is either zero, or else it is an element whose constant
term $\theta_0\in R$ lies in $\fm_R\!\setminus\!\fm^2_R$. The
assertion is obvious when $\theta=0$, and otherwise it follows from
claim \ref{cl_separ}. Next, by going down (corollary
\ref{cor_going-down-combin}), we see that :
$$
\hgt\,\fp A\geq\hgt\,\fp.
$$
On the other hand, notice that we have a natural identification
$\Spec\,\underline M_\xi\isom\Spec\,P$, especially
$\hgt\,\fp=\hgt\,\fq$, hence $\dim A/\fp A=\dim(A/\fq
A)^\wedge= \dim R+\dim(P\!\setminus\!\fq)-\eps$, where $\eps$ is
either $0$ or $1$ depending on whether $R$ does or does not contain
a field (corollary \ref{cor_new-dim}(ii)). Thus :
$$
\hgt\,\fp A\leq\dim A-\dim A/\fp A= \dim
P-\dim(P\!\setminus\!\fq)=\hgt\,\fp
$$
(corollary \ref{cor_consequent}(iii)), which completes the proof.

(ii): In this case, $\fq$ is a radical ideal of $P$, so it can
be written as a finite intersection of prime ideals of $P$
(lemmata \ref{lem_radical} and \ref{lem_face}(iii)); then the
assertion follows from (i) and lemma \ref{lem_intersect-ideals}
(details left to the reader).
\end{proof}

\begin{lemma}\label{lem_reduce-to-et}
Let $X$ be a scheme, $\underline M$ a fs log structure
on $X_\Zar$, and $\xi$ a geometric point of $X$. Then
$(X,\underline M)$ is regular at the point $|\xi|$ if
and only if\/ $\tilde u{}^*_X(X,\underline M)$ is regular
at $\xi$.
\end{lemma}
\begin{proof} Set $(Y,\underline N):=\tilde u{}^*_X(X,\underline M)$
(notation of \eqref{subsec_choose-a-top} : of course $Y=X$, but
the sheaf $\cO_Y$ is defined on the site $X_\et$, hence
$B:=\cO_{Y,\xi}$ is the strict henselization of
$A:=\cO_{\!X,|\xi|}$), and let $\alpha:\underline M{}_{|\xi|}\to B$
be the induced morphism of monoids; since $\alpha$ is local, it is
easily seen that $\underline N{}_\xi$ is isomorphic to the push-out
of the diagram
$\underline M{}_{|\xi|}\leftarrow\underline M{}^\times_{|\xi|}\to
B^\times$, especially $\dim\underline M{}_{|\xi|}=\dim\underline N{}_\xi$.
Moreover, $I(\xi,\underline N)=I(|\xi|,\underline M)B$, so that
$B_0:=B/I(\xi,\underline N)$ is the strict henselization of
$A_0:=A/I(|\xi|,\underline M)$ (\cite[Ch.IV, Prop.18.6.8]{EGA4});
hence $A_0$ is regular if and only if the same holds for $B_0$
(\cite[Ch.IV, Cor.18.8.13]{EGA4}). Finally, $\dim A=\dim B$ and
$\dim A_0=\dim B_0$ (\cite[Th.15.1]{Mat}). The lemma follows.
\end{proof}

\sset\subsubsection{}\label{subec_prepare-for-Kummer}
Let $(Y,\underline N)$ be a log scheme, $\bar y$ a $\tau$-point
of $Y$, and $\alpha:Q^{\log}_Y\to\underline N$ a fine and saturated
chart which is sharp at $\bar y$. Let also $\phi:Q\to P$ be an
injective morphism of monoids, with $P$ fine and saturated, and such
that $P^\times$ is a torsion-free abelian group. Define
$(X,\underline M)$ as the fibre product in the cartesian diagram of
log schemes :
\set\begin{equation}\label{eq_prepare-for-Kummer}
{\diagram
(X,\underline M) \ar[rr]^-f \ar[d] & &
(Y,\underline N) \ar[d]^h \\
\Spec(\Z,P) \ar[rr]^-{\Spec(\Z,\phi)} & & \Spec(\Z,Q)
\enddiagram}
\end{equation}
where $h$ is induced by $\alpha$; especially, $h$ is strict.

\begin{lemma}\label{lem_prepare-for-Kummer}
In the situation of \eqref{subec_prepare-for-Kummer},
let $\bar x$ be any $\tau$-point of $X$ such that
$f(\bar x)=\bar y$, and suppose that $(Y,\underline N)$
is regular at $\bar y$. Then $(X,\underline M)$ is
regular at $\bar x$.
\end{lemma}
\begin{proof} To begin with, we show :

\begin{claim}\label{cl_crit-1}
The natural map :
$$
\cO_{\!X,\bar x}\derotimes_PP/I\to \cO_{\!X,\bar x}\otimes_PP/I
$$
is an isomorphism in $\sD^-(\cO_{\!X,\bar x}\,\Mod)$, for every
ideal $I\subset P$ (notation of \eqref{subsec_restrict-scalar-mon}).
\end{claim}
\begin{pfclaim} It suffices to show that this map is an
isomorphism in $\sD^-(\cO_{Y,\bar y}\,\Mod)$, and in the latter
category we have a commutative diagram :
\set\begin{equation}\label{eq_bottom-is-loc}
{\diagram
(\cO_{Y,\bar y}\derotimes_QP)\derotimes_{P}P/I \ar[d] \ar[rr] & &
\cO_{Y,\bar y}\derotimes_QP/I \ar[d] \\
(\cO_{Y,\bar y}\otimes_QP)\derotimes_PP/I \ar[r] & (\cO_{Y,\bar
y}\otimes_QP)\otimes_PP/I \ar[r]^-\sim & \cO_{Y,\bar y}\otimes_QP/I.
\enddiagram}
\end{equation}
Since $\cO_{\!X,\bar x}$ is a localization of $\cO_{Y,\bar
y}\otimes_QP$, we are reduced to showing that the bottom arrow of
\eqref{eq_bottom-is-loc} is an isomorphism. However, the top arrow
of \eqref{eq_bottom-is-loc} is always an isomorphism. Moreover, on
the one hand, since $(Y,\underline N)$ is regular at $\bar y$, the
ring $\cO_{\!Y,\bar y}$ is $Q$-flat (proposition
\ref{prop_second-crit}(ii)), and on the other hand, since $\phi$ is
injective, $P/I$ is an integral $Q$-module, so the two vertical
arrows are isomorphisms as well, and the claim follows.
\end{pfclaim}

\begin{claim}\label{cl_crit-2}
$\cO_{\!X,\bar x}/I(\bar x,\underline M)$ is a regular ring.
\end{claim}
\begin{pfclaim} Let $\beta:P\to A:=\cO_{\!X,\bar x}$ be the morphism
deduced from $h$, and set :
$$
S:=\beta^{-1}(\cO^\times_{\!X,\bar x}) \qquad I(\bar
x,P):=P\!\setminus\!S.
$$
Then the $\Z[P]$-algebra $A$ is a localization of the
$\Z[P]$-algebra
$$
B:=S^{-1}(\cO_{Y,\bar y}\otimes_QP)
$$
and it suffices to show that $B/I(\bar x,\underline M)$ is regular.

It is easily seen that $I(\bar x,\underline M)=I(\bar x,P)A$ and
$\phi(\fm_Q)\subset I(\bar x,P)$; on the other hand, since $\alpha$
is sharp at the point $\bar y$, we have $I(\bar y,\underline
N)=\fm_Q\cO_{\!Y,\bar y}$, and $Q\setminus\fm_Q=\{1\}$. Let $p$ be
the residue characteristic of $\cO_{Y,\bar y}$; there follow
isomorphisms of $\Z_{(p)}$-algebras :
$$
B/I(\bar x,\underline M)B\simeq S^{-1}\cO_{Y,\bar
y}\otimes_QP/I(\bar x,P)\simeq \cO_{Y,\bar y}/I(\bar y,\underline
N)\otimes_{\Z_{(p)}} \Z_{(p)}[S^\gp].
$$
By assumption, $\cO_{Y,\bar y}/I(\bar y,\underline N)$ is regular,
hence we are reduced to showing that $\Z_{(p)}[S^\gp]$ is a
smooth $\Z_{(p)}$-algebra (\cite[Ch.IV, Prop.17.5.8(iii)]{EGA4}).
However, under the current assumptions $P^\gp$ is a free abelian
group of finite rank, hence the same holds for $S^\gp$, and the
contention follows easily.
\end{pfclaim}

In light of proposition \ref{prop_fla-criterion-point}(i), claims
\ref{cl_crit-1} and \ref{cl_crit-2} assert that conditions
(a) and (b) of proposition \ref{prop_second-crit}(i) are satisfied
for the $\tau$-point $\bar x$ of $(X,\underline M)$, so the latter
is regular at $\bar x$, as stated.
\end{proof}

\begin{theorem}\label{th_smooth-preserve-reg}
Let $f:(X,\underline M)\to(Y,\underline N)$ be a smooth morphism of
locally noetherian fs log schemes, $\xi$ a $\tau$-point of $X$,
and suppose that $(Y,\underline N)$ is regular at the point
$f(\xi)$. Then $(X,\underline M)$ is regular at the point $\xi$.
\end{theorem}
\begin{proof} In case $\tau=\Zar$, lemma \ref{lem_reduce-to-et} and
corollary \ref{cor_undercover}(ii) reduce the assertion to the
corresponding one for $\tilde u{}^*f$. Hence, we may assume that
$\tau=\et$. Next, since the assertion is local on $X_\tau$, we may
assume that $\underline N$ admits a fine and saturated chart
$\alpha:Q^{\log}_Y\to\underline N$ which is sharp at $f(\xi)$
(corollary \ref{cor_sharpy}(i)); then, by corollary
\ref{cor_charact-smoothness}, we may further assume that there
exist an injective morphism of fine and saturated monoids
$\phi:Q\to P$, such that $P^\times$ is torsion-free, and a cartesian
diagram as in \eqref{eq_prepare-for-Kummer}. Then the assertion
follows from lemma \ref{lem_prepare-for-Kummer}.
\end{proof}

\begin{corollary}\label{cor_smooth-preserve-reg}
Let $f:(X,\underline M)\to(Y,\underline N)$ be a smooth
morphism of locally noetherian fine log schemes. Then,
for every $y\in(Y,\underline N)_\tr$, the log scheme
$\Spec\,\kappa(y)\times_Y(X,\underline M)$ is geometrically
regular.
\end{corollary}
\begin{proof} The trivial log structure on $\Spec\,\kappa(y)$
is obviously saturated, so the same holds for the log
structure of $\Spec\,\kappa(y)\times_Y(X,\underline M)$.
Hence, the assertion is an immediate consequence of
theorem \ref{th_smooth-preserve-reg}.
\end{proof}

In the same vein we have :

\begin{proposition}\label{prop_we-got-a-situation}
Let $f:(X,\underline M)\to(Y,\underline N)$ be a
strict morphism of locally noetherian fs log schemes,
$\xi$ a $\tau$-point of $X$, and suppose that the
morphism of schemes $f_\xi:X(\xi)\to Y(f(\xi))$ is
flat. The following holds :
\begin{enumerate}
\item
If $(X,\underline M)$ is regular at the point $\xi$, then
$(Y,\underline N)$ is regular at the point $f(\xi)$.
\item
Conversely, if $(Y,\underline N)$ is regular at the point
$f(\xi)$, and the fibre $f^{-1}_\xi(f(\xi))$ is a regular
scheme (notation of \eqref{subsec_strict-loc-of-schs}),
then $(X,\underline M)$ is regular at the point $\xi$.
\end{enumerate}
\end{proposition}
\begin{proof} Set $A:=\cO_{Y,f(\xi)}$ and $B:=\cO_{\!X,\xi}$,
let $g:A\to B$ be the induced ring homomorphism, pick a
fine chart $\beta:P\to A$ for the log structure
$\underline N(f(\xi))$, and denote by $I\subset A$ the
ideal generated by the image of the maximal ideal of
$\underline N{}_{f(\xi)}$.

(i): By proposition \ref{prop_second-crit}(ii), the
local ring $B/IB$ is regular, $g\circ\beta:P\to B$ is
a fine chart for the log structure $\underline M(\xi)$,
and $B$ is $(g\circ\beta)$-flat. Since the induced
map $A/I\to B/IB$ is flat and local, also $A/I$ is
regular (\cite[Ch.0, Prop.17.3.3(i)]{EGAIV}), and
since $g$ is faithfully flat, $A$ is $\beta$-flat
(remark \ref{rem_when-is-it-flat}(iii)).
Then the assertion follows from proposition
\ref{prop_second-crit}(i).

(ii):  By proposition \ref{prop_second-crit}(ii), the
local ring $A/I$ is regular, and $A$ is $\beta$-flat.
Since $g$ is flat, the same holds for the induced map
$A/I\to B/IB$; it follows that also $B/IB$ is regular
(\cite[Ch.0, Prop.17.3.3(ii)]{EGAIV}), the induced map
$P\to B$ is a chart for $\underline M(\xi)$ and $B$ is
$(g\circ\beta)$-flat, whence the contention, again by
proposition \ref{prop_second-crit}(i).
\end{proof}

\begin{theorem}\label{th_reg-generizes}
Let $(X,\underline M)$ be a locally noetherian fs log scheme. Then
the subset $(X,\underline M)_\reg$ is closed under generization.
\end{theorem}
\begin{proof} Let $\xi$ be a $\tau$-point of $X$, with support
$x\in(X,\underline M)_\reg$, and let $\eta$ be a generization of
$\xi$, whose support is a strict generization $y$ of $x$. We have
to show that $(X,\underline M)$ is regular at the point $\eta$. Since
the assertion is local on $X$, we may assume that $(X,\underline M)$
admits a fine and saturated chart $\beta:P_X\to\underline M$, sharp
at the point $\xi$ (corollary \ref{cor_sharpy}(i)). Let
$\alpha:P\to\cO_{\!X,\eta}$ be the morphism deduced from $\beta_\eta$,
and set $\fp:=\alpha^{-1}\fm_\eta$, where $\fm_\eta\subset\cO_{\!X,\eta}$
is the maximal ideal. We consider the cartesian diagram of log schemes :
$$
\xymatrix{
(X',\underline M')_\circ \ar[r] \ar[d] &
(X,\underline M) \ar[d]^g \ar[d] \\
\Spec\La\Z,P/\fp\Ra \ar[r] & \Spec(\Z,P)
}$$
where $g$ is induced by $\beta$ (notation of \eqref{subsec_constant-log}).
Clearly $\xi$ and $\eta$ induce $\tau$-points on $X'$, which we
denote by the same names. We have natural identifications:
$$
\cO_{\!X,\xi}/I(\xi,\underline
M)\isom\cO_{\!X',\xi}/I(\xi,\underline M') \qquad \underline
M'_\eta=\cO^\times_{\!X',\eta} \qquad
\cO_{\!X,\eta}/I(\eta,\underline M)\isom\cO_{\!X',\eta}
$$
(\cite[Ch.IV, Prop.18.6.8]{EGA4}). Moreover, from proposition
\ref{prop_second-crit}(ii) we know that $\cO_{\!X,\xi}$ is
$P$-flat; then corollary \ref{cor_going-down-combin} yields
the inequality :
$$
\dim\cO_{\!X',\xi}=\dim\cO_{\!X,\xi}\otimes_PP/\fp\geq
\dim\cO_{\!X,\xi}-\hgt(\fp)=
\dim\cO_{\!X,\xi}/I(\xi,\underline M)+\dim P/\fp
$$
in other words : $\dim\cO_{\!X',\xi}=d(\xi,\underline M')$ (notation
of definition \ref{def_log-regular}), so $(X',\underline M')$ is
regular at $\xi$. By the same token, $\cO_{\!X,\eta}$ is $P_\fp$-flat;
from proposition \ref{prop_second-crit}(i), it follows that
$(X,\underline M)$ is regular at $\eta$ if and only if the same holds
for $(X',\underline M')$.

Hence we may replace $(X,\underline M)$ by $(X',\underline M')$,
and $P$ by $P\setminus\fp$, after which we may assume that
$y\in(X,\underline M)_\tr$. In this case :
\set\begin{equation}\label{eq_goes-to-invertible}
\alpha(P)\subset\cO_{\!X,\eta}^\times
\end{equation}
and we have to show that $\cO_{\!X,\eta}$ is a regular local ring,
or equivalently, that $\cO_{\!X,y}$ is regular
(\cite[Ch.IV, Cor.18.8.13(c)]{EGA4}). Now, if $P=0$, then $\cO_{\!X,x}$
is regular, and the assertion follows from \cite[Th.19.3]{Mat};
thus, we may assume that $P\neq 0$.

Denote by $Y$ the topological closure of $y$ in $X$, endowed with
its reduced subscheme structure. According to
\cite[Ch.II, Prop.7.1.7]{EGAII} (and
\cite[Ch.IV, Prop.18.8.8(iv)]{EGA4} in case $\tau=\et$), we may find
a local injective ring homomorphism $j:\cO_{Y,\xi}\to V$, where $V$
is a discrete valuation ring. Let also $\bar\beta:P\to\cO_{Y,\xi}$
be the morphism deduced from $\beta$.
Then \eqref{eq_goes-to-invertible} implies that
$j\circ\bar\beta(P)\subset V\!\setminus\!\{0\}$, hence
$j\circ\bar\beta$ extends to a homomorphism of groups
$P^\gp\to K^\times$, where $K^\times:=(V\!\setminus\!\{0\})^\gp$
is the multiplicative group of the field of fractions of $V$; after
composition with the valuation $K^\times\to\Z$ of $V$, there follows
a group homomorphism
$$
\phi:P^\gp\to\Z.
$$
Notice also that $j\circ\bar\beta$ is a local morphism, since the
same holds for $j$ and $\bar\beta$; consequently, $\phi(P^\gp)\neq\{0\}$.
Set $Q:=\phi^{-1}\N$; then $\phi(Q)$ is a non-trivial submonoid
of $\N$, and $\dim Q=\dim Q/\Ker\,\phi=\dim\phi(Q)=1$. Set
$$
T:=\Spec\,V
\qquad
(S,\underline Q):=\Spec(\Z,Q).
$$
Notice that $Q$ is saturated and fine (corollary
\ref{cor_fibres-are-fg}), so there exists an isomorphism :
$$
\Z^{\oplus r}\times\N\isom Q \qquad \text{for some $r\in\N$}
$$
(theorem \ref{th_structure-of-satu}(iii)); the latter determines a
chart :
\set\begin{equation}\label{eq_sharpy-charty}
\N_S\to\underline Q
\end{equation}
which is sharp at every $\tau$-point of $S$ localized
outside the trivial locus $\Spec(\Z,Q)_\tr$. We consider
the cartesian diagram of log schemes :
$$
\xymatrix{ (X',g^{\prime*}\underline Q) \ar[r]^-{g'}
\ar[d]_{f'} & \Spec(\Z,\underline Q) \ar[d]^f \\
(X,\underline M) \ar[r]^-g & \Spec(\Z,P) }
$$
where $f$ is the morphism of log schemes induced by the inclusion
map $\psi:P\to Q$. Notice that $f$ is a smooth morphism (proposition
\ref{prop_toric-smooth}); moreover, the restriction
$$
f_\tr:\Spec(\Z,Q)_\tr\to\Spec(\Z,P)_\tr
$$
of $f$, is just the morphism
$\Spec\,\Z[\psi^\gp]:\Spec\,\Z[Q^\gp]\to\Spec\,\Z[P^\gp]$, hence
it is an isomorphism of schemes. It follows that $f'$ is a smooth
morphism (proposition \ref{prop_sorite-smooth}(ii)), and its
restriction $f'_\tr$ to the trivial loci, is an isomorphism.

The homomorphism $j$ induces a morphism $h:T\to X$, such that
the closed point of $T$ maps to $x$ and the generic point maps
to $y$. By construction $g\circ h$ lifts to a morphism of schemes
$h':T\to S$, and the pair $(h,h')$ determines a morphism $T\to X'$.
Let $x',y'\in X'$ be the images of respectively the closed point
and the generic point of $T$, and choose $\tau$-points $\xi'$ and
$\eta'$ localized at $x'$ and respectively $y'$;
then the image of $x'$ in $X$ is the point $x$, therefore
$(X',g^{\prime*}\underline Q)$ is regular at the point $\xi'$,
by theorem \ref{th_smooth-preserve-reg}. Furthermore, $g'(\xi')$
lies outside $\Spec(\Z,Q)_\tr$, hence \eqref{eq_sharpy-charty}
induces a chart $\N_{X'}\to g^{\prime *}\underline Q$, which
is sharp at $\xi'$. In light of corollary \ref{cor_more-precisely},
we deduce that $\cO_{\!X',x'}$ is a regular ring, and then the
same holds also for $\cO_{\!X',y'}$ (\cite[Ch.0, Cor.17.3.2]{EGAIV}).
However, $y'$ lies in the trivial locus of
$(X',g^{\prime*}\underline Q)$, and its image in $X$ is $y$,
so by the foregoing the natural map $\cO_{\!X,y}\to\cO_{\!X',y'}$
is an isomorphism, and the contention follows.
\end{proof}

\begin{remark}\label{rem_correction-of-EGAII}
Notice that the proof of \cite[Chap.II, Prop.7.1.7]{EGAII} (that is
invoked in the proof of theorem \ref{th_reg-generizes}) is slightly
incorrect : indeed, with the notation of the proof of {\em loc.cit.},
it is implicitly assumed that $B/\fm B\neq 0$, which may fail,
{\em e.g.} take $A:=k[[x_2]]$ (for any field $k$), whose maximal
ideal is generated by $x_1:=x_2^2$ and $x_2$, so that $B=k((x_2))$
in this case. The proof can be amended by considering the scheme
$X:=\Spec\,A$ and the blow-up $\pi:E\to X$ of the coherent ideal
$\fm\cO_X$; since $\pi$ is proper and surjective, we see that there
exists $i\in\{1,\dots,n\}$ such the image of the restriction
$\Spec\,A[x_1/x_i,\dots,x_n/x_i]\to X$ of $\pi$ contains $\fm$, and
the argument then carries through with $B:=A[x_1/x_i,\dots,x_n/x_i]$.
See also \cite[Th.6.4.3]{Hu-Sw}.
\end{remark}

\sset\subsubsection{}\label{subsec_log-stratif}
Let $(X,\underline M)$ be any log scheme, $\xi$ a $\tau$-point
of $X$. There follows a continuous map
$$
\psi_\xi:X(\xi)\to T_\xi:=\Spec\,\underline M{}_\xi
$$
that sends the closed point of $X(\xi)$ to the closed point
of $T_\xi$. For every point $\fp\in T_\xi$, let $\bar{\{\fp\}}$
be the topological closure of $\{\fp\}$ in $T_\xi$; then
$$
X(\fp):=\psi_\xi^{-1}\bar{\{\fp\}}
$$
is a closed subset of $X(\xi)$, which we endow with its reduced
subscheme structure, and we set
$$
(X(\fp),\underline M(\fp)):=
((X(\xi),\underline M(\xi))\times_XX(\fp))_\red
$$
(notation of example \ref{ex_push-trivial}(iv)).
Notice that $U_\fp:=\psi_\xi^{-1}(\fp)$ is an open subset of
$X(\fp)$, for every $\fp\in T_\xi$. We call the family
$$
(U_\fp~|~\fp\in T_\xi)
$$
of
locally closed subschemes of $X(\xi)$, the
{\em logarithmic stratification\/} of $(X(\xi),\underline M(\xi))$.
For instance, $U_\emptyset=(X(\bar x),\underline M(\bar x))_\tr$.
More generally, it is clear from the definition that
\set\begin{equation}\label{eq_log-strat-tr}
U_\fp=(X(\fp),\underline M(\fp))_\tr
\qquad
\text{for every $\fp\in T_\xi$}.
\end{equation}

\begin{corollary}\label{cor_logar-strata}
In the situation of \eqref{subsec_log-stratif}, suppose that
$(X,\underline M)$ is regular at $\xi$. Then :
\begin{enumerate}
\item
The log scheme $(X(\fp),\underline M(\fp))$ is pointed regular,
for every $\fp\in T_\xi$.
\item
The scheme $X(\fp)$ is irreducible, and its codimension in $X$
equals the height of\/ $\fp$ in $T_\xi$, for every $\fp\in T_\xi$.
\item
The scheme $U_\fp$ is regular and irreducible, for every $\fp\in T_\xi$.
\end{enumerate}
\end{corollary}
\begin{proof} (i): By theorem \ref{th_reg-generizes}, it suffices
to show that $(X(\fp),\underline M(\fp))$ is pointed regular at
$\xi$, for every $\fp\in T_\xi$. However, say that
$X(\xi)=\Spec\,A$, and let $P\to\underline M(\xi)$ be a fine and
saturated chart, sharp at the $\tau$-point $\xi$. Then
$X(\fp)=\Spec\,A_0$, where $A_0:=A/\fp A$, and
$\underline M(\fp)=\underline N{}_\circ$, where $\underline N$
is the fs log structure deduced from the induced map
$\beta:P\!\setminus\!\fp\to A_0$. By proposition
\ref{prop_second-crit}(ii), the ring $A$ is $P$-flat;
then $A_0$ is $\beta$-flat, and the assertion follows from
proposition \ref{prop_second-crit}(i).

Next, (ii) is a rephrasing of corollary \ref{cor_same-height}(i),
and (iii) follows from (i),(ii) and \eqref{eq_log-strat-tr}, by
virtue of corollary \ref{cor_more-precisely}.
\end{proof}

\begin{proposition}\label{prop_log-struct-is-fixed}
Let $(X,\beta:\underline M\to\cO_{\!X})$ be a regular log scheme,
set $U:=(X,\underline M)_\tr$, and denote by $j:U\to X$ the open
immersion. Then the morphism $\beta$ induces identifications:
$$
\underline M\isom j_*\cO^\times_{\!U}\cap\cO_{\!X}
\qquad
\underline M^\gp\isom j_*\cO^\times_{\!U}.
$$
\end{proposition}
\begin{proof} Notice that the scheme $X$ is normal (corollary
\ref{cor_normal-and-CM}), hence both $j_*\cO^\times_{\!U}$
and $\cO_{\!X}$ are subsheaves of the sheaf $i_*\cO_{\!X_0}$,
where $X_0$ is the subscheme of maximal points of $X$, and
$i:X_0\to X$ is the natural morphism; so we may intersect
these two sheaves inside the latter.

In view of lemma \ref{lem_subsheaf-if-flat} and proposition
\ref{prop_second-crit}(ii), we know already that $\beta$ is injective,
and clearly the image of $\beta$ lands in $j_*\cO^\times_{\!U}$,
so it remains only to show that $\beta$ (resp. $\beta^\gp$) induces
an epimorphism onto $j_*\cO^\times_{\!U}\cap\cO_{\!X}$ (resp.
onto $j_*\cO^\times_{\!U}$). The assertions can be checked
on the stalks, hence let $\xi$ be any $\tau$-point of $X$; to
begin with, we show :

\begin{claim}\label{cl_reflex-from-triv}
The induced map
$\beta_\xi^\gp:\underline M^\gp_\xi\to(j_*\cO^\times_{\!U})_\xi$
is a surjection.
\end{claim}
\begin{pfclaim} Set $A:=\cO_{\!X,\xi}$; notice that
$(X(\xi),\underline M(\xi))_\tr$ is the complement of the union
of the finitely many closed subsets of the form
$Z_\fp:=\Spec\,A/\fp A$, where $\fp\subset\underline M_\xi$ runs
over the prime ideals of height one. Due to corollary
\ref{cor_same-height}(i), each $Z_\fp$ is an irreducible divisor
in $\Spec\,A$. Now, let $s\in(j_*\cO^\times_{\!U})_\xi$; then the
divisor of $s$ is of the form $\sum_{\hgt\fp=1}n_\fp[\fp]$, for
some $n_\fp\in\Z$. By lemmata \ref{lem_simple-charts-top}(ii)
and \ref{lem_simple-charts}(i), $\underline M{}^\sharp_\xi$
is a fine and saturated monoid; then lemma \ref{lem_rflx_and_quot}
and proposition \ref{prop_classify-reflex-mon} say that the
fractional ideal $I:=\bigcap_{\hgt\,\fp=1}\fm_\fp^{n_\fp}$ of
$\underline M{}_\xi$ is reflexive, and therefore the fractional ideal
$IA$ of $A$ is reflexive as well (lemma \ref{lem_rflx-rflx}(ii)).
Since $A$ is normal, by considering the localizations of
$I\cO_{\!X,\xi}$ at the prime ideals of height one of $A$, we deduce
easily that $IA=sA$ (\cite[Ch.VII, \S4, n.2, Cor. du Th.2]{BouAC}).
Thus, we may write $s=\sum_{i=1}^n\beta^\gp_\xi(x_i)a_i$ for certain
$a_i\in A$ and $x_i\in I$ (for $i=1,\dots,n$). Then we must have
$\beta^\gp_\xi(x_i)a_i\notin s\cdot\fm_\xi$ for at least one index
$i\leq n$ (where $\fm_\xi\subset A$ is the maximal ideal); for such
$i$, it follows that $s^{-1}\cdot\beta^\gp_\xi(x_i)\in A^\times$,
whence $s\in\underline M{}^\gp_\xi$, as required.
\end{pfclaim}

Now, let $s\in(j_*\cO^\times_{\!U})_\xi\cap A$; by claim
\ref{cl_reflex-from-triv}, we may find $x\in\underline M^\gp_\xi$
such that $\beta_\xi^\gp(x)=s$. To conclude, it suffices to
show that $x\in\underline M{}_\xi$. By theorem
\ref{th_structure-of-satu}(i), we are reduced to showing that
$x\in(\underline M{}_\xi)_\fp$ for every prime ideal
$\fp\subset\underline M{}_\xi$ of height one. Set $\fq:=\fp A$;
then $\fq$ is a prime ideal of height one (corollary
\ref{cor_same-height}(i)), and $\beta_\xi$ extends to a well defined
map of monoids $\beta_\fp:(\underline M{}_\xi)_\fp\to A_\fq$.
Furthermore, $A_\fq$ is a discrete valuation ring (corollary
\ref{cor_normal-and-CM}), whose valuation we denote
$v:A_\fq\setminus\{0\}\to\N$.
In view of theorem \ref{th_structure-of-satu}(ii), we see
that $x\in(\underline M{}_\xi)_\fp$ if and only if
$(v\circ\beta_\fp)^\gp(x)\geq 0$. However, clearly $v(s)\geq 0$,
whence the contention.
\end{proof}

\begin{remark}\label{rem_fully-regular}
Let $\mathbf{reg.log}_\tau$ denote the full subcategory of
$\blog_\tau$ whose objects are the regular log schemes.
As an immediate consequence of proposition
\ref{prop_log-struct-is-fixed} (and of remark
\ref{rem_trivial-is-strict}(i)) we see that the forgetful
functor $F$ of \eqref{subsec_special-schs} restricts to a
fully faithful functor
$$
\mathbf{reg.log}_\tau\to\mathbf{Open}
\qquad
(X,\underline M)\mapsto((X,\underline M)_\tr\to X)
$$
where $\mathbf{Open}$ is full subcategory of $\sMorph(\Sch)$
whose objects are the open immersions.
\end{remark}

\sset\subsubsection{}\label{subsec_dualzing-logreg}
Let $(X,\underline M)$ be a quasi-compact and regular log
scheme. Set $Z:=X\setminus(X,\underline M)_\tr$; then $Z$ is
a closed subset of $X$ (lemma \ref{lem_simple-charts}(ii)),
and we endow it with its reduced subscheme structure.
Denote by $\omega_X\subset\cO_{\!X}$ the sheaf of ideals
corresponding to $Z$. We have :

\begin{theorem}\label{th_log-reg-dualized}
In the situation of \eqref{subsec_dualzing-logreg}, the
complex $\omega_X[0]$ is dualizing on $X$. (See definition
{\em\ref{def_dual-complex}}.)
\end{theorem}
\begin{proof} In light of proposition \ref{prop_fishy}, we
may assume that $X$ is local, say $X=\Spec\,A$, for a local
noetherian ring $A$. Next, pick any $\tau$-point $\xi$
localized at the closed point of $X$; by virtue of proposition
\ref{prop_Goren-pullback} and
\cite[Ch.IV, Prop.18.8.8(iv) and Prop.18.8.12(ii)]{EGA4},
we may further reduce to the case where
$(X,\underline M)=(X(\xi),\underline M(\xi))$. Let
$X^\wedge:=\Spec\,A^\wedge$, where $A^\wedge$ denotes the
completion of the local ring $A$, and denote by $f:X^\wedge\to X$
the natural morphism; by applying again proposition
\ref{prop_Goren-pullback}, we reduce to checking that
$f^*\omega_X[0]$ is dualizing on $X^\wedge$. However,
theorem \ref{th_charact-log-regular} says that there
exist a local ring homomorphism $R\to A^\wedge$, with
$R$ a complete regular local ring, and a fine and saturated
chart $P_X\to\underline M$, sharp at the closed point of
$X$, such that the induced continuous ring homomorphism
$\phi:R[[P]]\to A^\wedge$ is an isomorphism if $A$ contains a
field, and otherwise it is a surjection, whose kernel is
generated by a regular element $\theta\in R[[P]]$. With this
notation, a simple inspection shows that $Z$ is the union
of the closed subsets $X(\fp)\subset X$, where $\fp\subset P$
ranges over all prime ideals $\fp\neq\emptyset$ (notation
of \eqref{subsec_log-stratif}).
Then, corollary \ref{cor_same-height}(i) says that
$\omega_X$ is the intersection of the ideals of $\cO_{\!X}$
of the form $\fp\cO_{\!X}$, where $\fp\subset P$ is an
arbitrary non-empty prime ideal. Denote by $P^\circ\subset P$
the intersection of all the (finitely many) non-empty prime
ideals of $P$; lemma \ref{lem_intersect-ideals} and
proposition \ref{prop_second-crit}(ii) imply that
$\omega_X=P^\circ\cdot\cO_{\!X}$ (details left to the reader),
and since $f$ is a flat morphism, we deduce that
$f^*\omega_X=P^\circ\cdot\cO_{\!X^\wedge}$.
The ideal $P^\circ$ can also be described as follows.
Set $V:=P^\gp_\R$, denote by $\sigma$ the unique strictly
convex polyhedral cone such that $P=\sigma\cap P^\gp$, and
let $\sigma^\circ$ be the topological interior of $\sigma$
in $V$; then it is easily seen that $P^\circ=P\cap\sigma^\circ$.

Set $Y:=\Spec\,R[P]$, $\omega_Y:=P^\circ\cdot\cO_Y$,
and let $y\in Y$ be the image of the closed point of
$Y^\wedge:=\Spec\,R[[P]]$, under the natural morphism
$g:Y^\wedge\to Y$. Theorem \ref{th_Stanley}(ii) says that
$\omega_Y[0]$ is dualizing on $Y$; then, since $g$ is flat,
propositions \ref{prop_Goren-pullback} and \ref{prop_fishy}
imply that $g^*\omega_Y[0]$ is dualizing on $Y^\wedge$.

Now, suppose that $A$ contains a field; then $\phi$ induces
an isomorphism $X^\wedge\isom Y^\wedge$ which identifies
$f^*\omega_X$ with $g^*\omega_Y$, whence the contention.

Lastly, if $A$ does not contain a field, $\phi$ induces
a regular closed immersion $i:X^\wedge\to Y^\wedge$, so
$i^!(g^*\omega_Y[0])$ is dualizing on $X^\wedge$ (lemma
\ref{lem_transit-dual}(i)). The latter is the complex
of $\cO_{\!X^\wedge}$-modules arising from the complex
of $A^\wedge$-modules
$L^\bullet:=R\Hom^\bullet_{R[[P]]}(A^\wedge[0],R[[P^\circ]])$
(corollary \ref{cor_Ext-loc=glob}(ii); here
$R[[P^\circ]]:=P^\circ\cdot R[[P]]$). However, $A^\wedge[0]$
is naturally isomorphic (in $\sD(R[[P]]\Mod)$) to the
Koszul complex $\bK_\bullet(\theta)$ (notation of
remark \ref{rem_koszul-alg}(ii) : this is a bounded
complex of free $R[[P]]$-modules), so we see that
$$
L^\bullet\simeq(R[[P^\circ]]\otimes_{R[[P]]}A^\wedge)[-1].
$$
However, we have natural $A^\wedge$-linear isomorphisms :
$$
\begin{aligned}
R[[P^\circ]]\otimes_{R[[P]]}A^\wedge\isom\: &
R[P^\circ]\otimes_{R[P]}A^\wedge & & \text{(since $g$ is flat)} \\
\isom\: & \Z[P^\circ]\otimes_{\Z[P]}A\otimes_AA^\wedge \\
\isom\: & (P^\circ\cdot A)\otimes_AA^\wedge & &
\text{(by proposition \ref{prop_second-crit}(ii))}
\end{aligned}
$$
which shows that $i^!(g^*\omega_Y[0])\simeq f^*\omega_X[-1]$,
and concludes the proof of the theorem.
\end{proof}

\subsection{Resolution of singularities of regular log schemes}
Most of this section concerns results that are special
to the class of log schemes over the Zariski topology; these
are then applied to \'etale fs log structures, after we have
shown that every such log structure admits a logarithmic blow
up which descends to the Zariski topology (proposition
\ref{subsec_descends-to-Zar}). The same statement -- with
an unnecessary restriction to regular log schemes -- can be
found in the article \cite{Niz} by W.Niziol : see theorem 5.6
of {\em loc.cit.} We mainly follow her treatment, except for
fleshing out some details, and correcting some inaccuracies.

Therefore, {\em we let here $\tau=\Zar$, and all log structures
considered in this section until \eqref{subsec_down-from-log}
are defined on the Zariski sites of their underlying schemes}.

\sset\subsubsection{}\label{subsec_resolve-logscheme}
Let $f:(Y,\underline N)\to(X,\underline M)$ be any morphism of
log schemes; we remark that $f$ is a morphism of monoidal spaces,
{\em i.e.} for every $y\in Y$, the map
$\underline M{}_{f(y)}\to\underline N{}_y$ induced by $\log f$,
is local. Indeed, it has been remarked in
\eqref{subsec_up_and_down-log} that the natural map
$\underline M{}_{f(y)}\to(f^*\underline M)_y$ is local, and on
the other hand, a section $s\in(f^*\underline M)_y$ is invertible
if and only if its image in $\cO_{Y,y}$ is invertible, if and
only if $\log f_y(s)$ is invertible in $\underline N{}_y$, whence
the contention. We shall denote
$$
f^\sharp:(Y,\underline N)^\sharp\to(X,\underline M)^\sharp
$$
the morphism of sharp monoidal spaces induced by $f$ in the
obvious way ({\em i.e.} the underlying continuous map is the
same as the continuous map underlying $f$, and
$\log f^\sharp:f^*\underline M^\sharp\to\underline N^\sharp$
is $(\log f)^\sharp$ : see definition \ref{def_mon-spaces}).

\sset\subsubsection{}\label{subsec_category-K}
We let $\cK$ be the category whose objects are all data of the form
$\underline X:=((X,\underline M),F,\psi)$, where $(X,\underline M)$
is a log scheme, $F$ is a fan, and $\psi:(X,\underline M)^\sharp\to F$
is a morphism of sharp monoidal spaces, such that
$\log\psi:\psi^*\cO_F\to\underline M^\sharp$ is an isomorphism.
The morphisms :
$$
((Y,\underline N),F',\psi')\to\underline X
$$
in $\cK$ are all the pairs $(f,\phi)$, where
$f:(Y,\underline N)\to(X,\underline M)$ is a morphism of log schemes,
and $\phi:F'\to F$ is a morphism of fans, such that the diagram
$$
\xymatrix{ (Y,\underline N)^\sharp \ar[r]^-{f^\sharp}
\ar[d]_{\psi'} & (X,\underline M)^\sharp \ar[d]^\psi \\
F' \ar[r]^-\phi & F}
$$
commutes. Especially, notice that lemma \ref{lem_check-iso}
implies the identity :
\set\begin{equation}\label{eq_strict-identity}
\mathrm{Str}(f)=\psi^{\prime-1}\mathrm{Str}(\phi).
\end{equation}
(notation of definitions \ref{def_mon-spaces}(ii) and
\ref{def_trivial-locus}(ii)).  We shall say that an object
$\underline X$ is {\em locally noetherian}, if the same holds
for the scheme $X$. Likewise, a morphism $(f,\phi)$ in $\cK$
is {\em quasi-compact}, (resp. {\em quasi-separated}, resp.
{\em separated}, resp. {\em locally of finite type}, resp.
{\em of finite type}) if the same holds for the morphism of
schemes underlying $f$. We say that $f$ is {\em \'etale}, if
the same holds for the morphism of log schemes underlying $f$.
Also, the {\em trivial locus\/} of $\underline X$ is defined
as the subset $(X,\underline M)_\tr$ of $X$.

\sset\subsubsection{}\label{subsec_base-change-in-cK}
There is an obvious (forgetful) functor :
$$
F:\cK\to\Sch
\qquad
((X,\underline M),F,\psi)\mapsto X
$$
which is a fibration. More precisely, every morphism of schemes
$S'\to S$ induces a {\em base change\/} functor (notation of
\eqref{subsec_fibreovercat}) :
$$
F\cK/S\to F\cK/S'
\qquad
\underline X\mapsto S'\times_S\underline X
$$
unique up to natural isomorphism of functors. Namely, for
$\underline X:=((X,\underline M),F,\psi)$ one lets
$$
S'\times_S\underline X:=(S'\times_S(X,\underline M),F,\psi')
$$
where $\psi':=\psi\circ\pi^\sharp$, and
$\pi:S'\times_S(X,\underline M)\to(X,\underline M)$ is the natural
projection.

\begin{example}\label{ex_fan-andlogscheme}
(i)\ \  
Let $P$ be a monoid, $R$ a ring, and set
$(S,P^{\log}_{\!S}):=\Spec(R,P)$ (see \eqref{subsec_constant-log}).
The unit of adjunction $\eps_P:P\to R[P]$ determines
a unique morphism of sharp monoidal spaces
$$
\psi_P:\Spec(R,P)^\sharp\to T_P:=(\Spec\,P)^\sharp
$$
(proposition \ref{prop_Spec-fan}), and we claim that
$\underline S:=(\Spec(R,P),T_P,\psi_P)$ is an object of $\cK$.

The assertion can be checked on the stalks, hence let $\xi$
be a point of $S$; then $\xi$ corresponds to a prime ideal
$\fp_\xi\subset R[P]$, and by inspecting the definitions,
we see that
$\psi_P(\xi)=\eps_P^{-1}(\fp_\xi)=\fq_\xi:=P\cap\fp_\xi\in\Spec\,P$.
Again, a direct inspection shows that the morphism
$$
(\log\psi_P)_\xi:\cO_{T_P,\fq_\xi}\to(P^{\log}_{\!S})^\sharp_\xi
$$
is none else than the natural isomorphism
$P^\sharp_{\fq_\xi}\isom P^{\log}_{\!S,\xi}/\cO^\times_{\!S,\xi}$
deduced from $\eps_P$ and the natural identifications
$$
P^{\log}_{\!S,\xi}=\cO^\times_{\!S,\xi}\otimes_{P\!\setminus\!\fq_\xi}P
=\cO^\times_{\!S,\xi}\otimes_{P_{\fq_\xi}^\times}P_{\fq_\xi}.
$$

(ii)\ \ 
The construction of $\underline S$ is clearly functorial in $P$.
Namely, say that $\lambda:P\to Q$ is a morphism of monoids, and
set $S':=\Spec\,R[Q]$, $T_Q:=(\Spec\,Q)^\sharp$. There follows
a morphism
\set\begin{equation}\label{eq_clearly-fun}
\underline S':=(\Spec(R,Q),T_Q,\psi_Q)\to\underline S
\end{equation}
in $\cK$, whose underlying morphism of log schemes is
$\Spec(R,\lambda)$, and whose underlying morphism of fans
is just $(\Spec\,\lambda)^\sharp$.
\end{example}

\sset\subsubsection{}\label{subsec_fan-of-a-logsch}
Example \ref{ex_fan-andlogscheme} can be globalized to more general
log schemes, at least under some additional assumptions. Namely,
let $(X,\underline M)$ be a regular log scheme. For every point
$x$ of $X$, let $\fm_x\subset\cO_{\!X,x}$ be the maximal ideal;
we set :
$$
F(X):=\{x\in X~|~ I(x,\underline M)=\fm_x\}
$$
(notation of \eqref{subsec_define-new-dim}), and we endow $F(X)$
with the topology induced from $X$. The {\em fan of
$(X,\underline M)$} is the sharp monoidal space
$$
F(X,\underline M):=(F(X),\underline M_{|F(X)})^\sharp.
$$
We wish to show that $F(X,\underline M)$ is indeed a fan. Let
$U\subset X$ be any open subset; to begin with, it is clear that
$F(U,\underline M_{|U})$ is naturally an open subset of
$F(X,\underline M)$; hence the contention is local on $X$,
so we may assume that $X$ is affine, say $X=\Spec\,A$ for a
noetherian ring $A$, and that $\underline M$ admits a finite
chart $P_X\to\underline M$; denote by $\beta:P\to A$ the induced
morphism of monoids. Let now $\xi$ be a point of $X$, and
$\fp_\xi\subset A$ the corresponding prime ideal; set
$\fq_\xi:=\beta^{-1}\fp_\xi\in\Spec\,P$.
By inspecting the definitions, it is easily seen that
$I(\xi,\underline M)=\fq_\xi A_{\fp_\xi}$. Therefore, for any
prime ideal $\fq\subset P$, let $V(\fq)_{\max}$ be the finite
set consisting of all the maximal points of the closed subset
$V(\fq):=\Spec\,A/\fq A$ ({\em i.e.} the minimal prime ideals
of $A/\fq A$); in light of corollary \ref{cor_same-height}(i),
it follows easily that
\set\begin{equation}\label{eq_describe-fan-of-X}
F(X)=\bigcup_{\fq\in\Spec\,P}V(\fq)_{\max}
\end{equation}
especially, $F(X)$ is a finite set (lemma \ref{lem_face}(iii)).
Let $t\in F(X)$ be any element, and denote by $U(t)\subset F(X)$
the subset of all generizations of $t$ in $F(X)$; as a corollary,
we see that $U(t)$ is an open subset of $F(X)$. Moreover,
\eqref{eq_describe-fan-of-X} also implies that :
$$
U(t)=\bigcup_{\fq\in\Spec\,P}
(\Spec\,\cO_{\!X,t}/\fq\cO_{\!X,t})_{\max}.
$$
However, corollary \ref{cor_same-height}(i) says that
$\Spec\,\cO_{\!X,t}/\fq\cO_{\!X,t}$ is irreducible for every
$\fq\in\Spec\,P$, therefore the set $U(t)$ is naturally identified
with a subset of $\Spec\,P$. Furthermore, arguing as in example
\ref{ex_fan-andlogscheme}(i) we find a natural isomorphism
$P^\sharp_{\!\fq_t}\isom\cO_{\!F(X,\underline M),t}$.
Moreover, if $t'\in U(t)$ is any other point, the specialization
map $\cO_{\!F(X,\underline M),t}\to\cO_{\!F(X,\underline M),t'}$
corresponds -- under the above isomorphism -- to the natural
morphism $P^\sharp_{\!\fq_t}\to P^\sharp_{\!\fq_t'}$ induced
by the localization map $P_{\!\fq_t}\to P_{\!\fq_t'}$. This
shows that the open monoidal subspace $(U(t),\underline M_{|U(t)})$
is naturally isomorphic to $(\Spec\,P_{\!\fq_t})^\sharp$, hence
$F(X,\underline M)$ is a fan, as stated.

\begin{remark}\label{rem_discuss-fan-logreg}
(i)\ \
The discussion in \eqref{subsec_fan-of-a-logsch} shows more
precisely that if :
\begin{enumerate}
\alphaenu
\item
$(U,\underline M)$ is a log scheme with $U=\Spec\,A$ affine,
\item
there exists a morphism $\beta:P\to A$ from a finitely generated
monoid $P$, inducing the log structure $\underline M$ on $U$, and
\item
$\beta(\fq)A$ is a prime ideal for every $\fq\in\Spec\,P$
\end{enumerate}
then $F(U,\underline M)$ is naturally identified with
$(\Spec\,P)^\sharp$. In this situation, denote by
$$
f_\beta:(U,\underline M)^\sharp\to T_P:=(\Spec\,P)^\sharp
$$
the morphism of sharp monoidal spaces deduced from $\beta$
(proposition \ref{prop_Spec-fan}); by inspecting the definitions,
we see that the associated map
$\log f_\beta:f^*_\beta\cO_{T_P}\to\underline M^\sharp$
is an isomorphism. Via the foregoing natural identification, there
results a morphism $\pi_U:(U,\underline M)^\sharp\to F(U,\underline M)$
which can be described without reference to $P$. Indeed, let $x\in U$
be any point, and $\fp\in\Spec\,A$ the corresponding prime ideal; by
inspecting the definitions we find that
$$
\pi_U(x)=I(x,\underline M)
\ \ 
\text{which is the largest prime
ideal in $\Spec\,A_\fp\cap F(U,\underline M)$}.
$$
Especially, $\pi_U(x)$ is a generization of $x$, and the inverse
of $(\log\pi_U)_x:\cO_{\!F(U,\underline M),\pi_U(x)}\isom
\underline M{}_x^\sharp$ is induced by the
specialization map $\underline M{}_x\to\underline M{}_{\pi_U(x)}$
(which induces an isomorphism on the associated sharp quotient
monoids).

(ii)\ \
Finally, it follows easily from corollary \ref{cor_same-height}(i)
and \cite[Ch.IV, Cor.8.4.3]{EGAIV-3}, that every regular log scheme
$(X,\underline M)$ admits an affine open covering
$X=\bigcup_{i\in I}U_i$, such that each $(U_i,\underline M{}_{|U_i})$
fulfills conditions (a)--(c) above, and the intrinsic description
in (i) shows that the morphisms $\pi_{U_i}$ glue to a well defined
morphism $\pi_X:(X,\underline M)^\sharp\to F(X,\underline M)$ of
sharp monoidal spaces, such that the datum
$$
\cK(X,\underline M):=((X,\underline M),F(X,\underline M),\pi_X)
$$
is an object of $\cK$.

(iii)\ \
Notice that $\pi_X^{-1}(t)$ is an irreducible locally closed
subset of $X$ of codimension equal to the height of $t$, for
every $t\in F(X,\underline M)$; also $\pi_X^{-1}(t)$ is a
regular scheme, for its reduced subscheme structure. Indeed,
we have already observed that $t$ is the unique maximal point of
$\pi_X^{-1}(t)$, and then the assertion follows immediately from
corollary \ref{cor_logar-strata}(ii,iii). Furthermore, the
inclusion map $j:F(X,\underline M)\to X$ is a continuous section
of $\pi_X$, and notice that $j\circ\pi_X(U)=j^{-1}U$ for every
open subset $U\subset X$, since $j$ maps every $t\in F(X,\underline M)$
to the unique maximal point of $\pi_X^{-1}(t)$. It follows easily
that $\pi_X$ is an open map.
\end{remark}

\sset\subsubsection{}\label{subsec_not_quite_fibration}
Denote by $\cK_\intg$ the full subcategory of $\cK$ whose objects
are the data $((X,\underline M),F,\psi)$ such that $\underline M$
is an integral log structure, and $F$ is an integral fan. There
is an obvious functor
$$
\cK_\intg\to\mathbf{int.Fan} \quad :\quad ((X,\underline
M),F,\psi)\mapsto F
$$
to the category of integral fans, which shall be used to construct
useful morphisms of log schemes, starting from given morphisms of
fans. This technique rests on the following three results :

\begin{lemma}\label{lem_soon-to-be-named}
Let $((X,\underline M),F,\psi)$ be an object of\/ $\cK_\intg$, with
$F$ locally fine. Then, for each point $x\in X$ there exist an open
neighborhood $U\subset X$ of $x$, a fine chart
$Q_U\to\underline M{}_{|U}$ (for some fine monoid $Q$ depending on
$x$), and an isomorphism of monoids $Q^\sharp\isom\cO_{F,\psi(x)}$.
\end{lemma}
\begin{proof} The assertion is local on $X$, hence we may
assume that $F=(\Spec\,P)^\sharp$ for a fine monoid $P$. In
this case, $\psi$ is determined by the corresponding map
$$
\bar\beta:P_X\to\underline M^\sharp.
$$
Indeed, for any $x\in X$, the point $\psi(x)$ is the prime ideal
$\bar\beta{}^{-1}_x(\fm_x)\subset P$, where
$\fm_x\subset\underline M^\sharp_x$ is the maximal ideal;
moreover :
$$
\cO_{\!F,\psi(x)}=P/\bar S_x
\qquad\text{where}\qquad
\bar S_x:=\bar\beta{}^{-1}_x(1)
$$
and -- under this identification -- the isomorphism
$\log\psi_x:\cO_{\!F,\psi(x)}\isom\underline M^\sharp_x$
is deduced from $\bar\beta_x$ in the obvious way.

Now, let $x\in X$ be any point; after replacing $X$ by the open
subset $\psi^{-1}U(\psi(x))$, we may assume that
$P=\cO_{\!F,\psi(x)}$ (notation of \eqref{subsec_height-in-T}), and
by assumption $\log\psi_x:P\to\underline M^\sharp_x$ is an
isomorphism. Pick a surjection $\alpha:\Z^{\oplus r}\to P^\gp$, and
let $Q$ be the pull-back in the cartesian diagram :
$$
\xymatrix{
Q \ar[r] \ar[d] &
\Z^{\oplus r} \ar[d]^{(\log\psi_x)^\gp\circ\alpha} \\
\underline M{}^\sharp_x \ar[r] &
\underline M{}^\gp_x/\underline M{}^\times_x.
}$$
By construction, $\alpha$ restricts to a morphism of monoids
$\theta:Q\to P$, inducing an isomorphism $Q^\sharp\isom P$,
whence an isomorphisms of fans
$(\Spec\,P)^\sharp\isom(\Spec\,Q)^\sharp$, and $Q$ is fine,
by corollary \ref{cor_fibres-are-fg}. Moreover, the choice of
a lifting $\Z^{\oplus r}\to\underline M{}^\gp_x$ of
$(\log\psi_x)^\gp\circ\alpha$ determines a morphism
$$
\beta_x:Q\to\underline M{}^\sharp_x
\times_{\underline M{}^\gp_x/\underline M{}^\times_x}
\underline M{}^\gp_x=\underline M{}_x
$$
which lifts $\log\psi_x\circ\theta$ (here it is needed that
$\underline M$ is integral). Next, after replacing $X$ by
an open neighborhood of $x$, we may assume that $\beta_x$ extends
to a map of pre-log structures $\beta:Q_X\to\underline M$ lifting
$\bar\beta$ (lemma \ref{lem_simple-charts-top}(iv.b),(v)).
It remains only to show that $\beta$ is a chart for $\underline M$,
which can be checked on the stalks. Thus, let $y\in X$ be any point,
and set $S_y:=\beta^{-1}_y\underline M_y^\times$; with the foregoing
notation, the stalk $Q^{\log}_{\!X,y}$ of the induced log structure
is naturally isomorphic to
$(S^{-1}_yQ\times\cO_{\!X,y}^\times)/S^\gp_y$, therefore
$$
(Q^{\log}_{\!X,y})^\sharp\simeq Q/S_y\simeq P/\bar S_y
$$
and the induced map $P/\bar S_y\to\underline M{}_y^\sharp$ is again
deduced from $\bar\beta_y$, so it is an isomorphism; then the
same holds for $\beta^{\log}_y:Q^{\log}_{\!X,y}\to\underline M{}_y$
(lemma \ref{lem_check-iso}).
\end{proof}

\begin{lemma}\label{lem_redux-more}
Let $\mu:P\to P'$ be a morphism of integral monoids such that
$\mu^\gp$ is surjective,
$\phi:(\Spec\,P')^\sharp\to(\Spec\,P)^\sharp$ the induced morphism
of affine fans, and denote by
$$
\lambda:P\to Q:=P^\gp\times_{P^{\prime\gp}}P'
$$
the map of monoids determined by $\mu$ and the unit of adjunction
$P\to P^\gp$. Then :
\begin{enumerate}
\item
The natural projection $Q\to P'$ induces an isomorphism
$\omega:(\Spec\,P')^\sharp\to(\Spec\,Q)^\sharp$ of affine fans,
such that $(\Spec\,\lambda)^\sharp\circ\omega=\phi$.
\item
The induced morphism \eqref{eq_clearly-fun} in\/ $\cK_{\intg}$
(with $R:=\Z$) is\/ $\mathbf{int.Fan}$-cartesian.
\item
If $P$ and $P'$ are fine, the morphism \eqref{eq_clearly-fun} is
\'etale.
\end{enumerate}
\end{lemma}
\begin{proof} (See \eqref{sec_phi-cartesian} for generalities
concerning inverse images and cartesian morphisms relative to a
functor.) Notice first that, since $\mu^\gp$ is surjective, the
projection $Q\to P'$ induces an isomorphism
$Q^\sharp\isom P^{\prime\sharp}$, whence (i).

(ii): Notice that the log structure of $\Spec(\Z,Q)$ is integral,
by lemma \ref{lem_simple-charts-top}(iii). Next, set
$$
F':=(\Spec\,P')^\sharp.
$$
Define $\underline S$, $\underline S'$ as in example
\ref{ex_fan-andlogscheme}(ii) (with $R:=\Z$), let
$g:((Y,\underline N),F'',\psi_Y)\to\underline S$ be any morphism
of $\cK_\intg$, and $\phi':F''\to F'$ a morphism of integral fans,
such that the image of $g$ in $\mathbf{int.Fan}$ equals
$\phi\circ\phi'$; we must show that $g$ factors through a unique
morphism $h:((Y,\underline N),F'',\psi_Y)\to\underline S'$, whose
image in $\mathbf{int.Fan}$ equals $\phi'$. As usual, we may reduce
to the case where $Y=\Spec\,A$ is affine, $F''=\Spec\,P''$ is an
affine fan for a sharp integral monoid $P''$, and $\psi_Y$ is given
by a map of sheaves $P''_Y\to\underline N^\sharp$. In such
situation, the morphism $(Y,\underline N)\to\Spec(\Z,P)$
underlying $g$ is determined by a morphism of monoids
$P\to\underline N(Y)$, or which is the same, a map of sheaves
$\gamma_Y:P_{\!Y}\to\underline N$; likewise, $\phi'$ is given by a
morphism of monoids $P'\to P''$, and composing with $\psi_Y$, we get
a map of sheaves $\alpha_Y:P'_{\!Y}\to\underline N^\sharp$. Finally,
the condition that $g$ lies over $\phi\circ\phi'$ translates as the
commutativity of the following diagram of sheaves :
$$
\xymatrix{ P^\gp_{\!Y} \ar[r] \ar[d]_{\gamma_Y^\gp} &
P^{\prime\gp}_{\!Y} \ar[d] &
\ar[l] P'_{\!Y} \ar[d]^{\alpha_Y} \\
\underline N^\gp \ar[r] & \underline N^\gp/\underline N^\times &
\ar[l] \underline N^\sharp. }
$$
There follows a unique morphism of sheaves :
\set\begin{equation}\label{eq_from-Q}
Q_Y\to\underline N^\gp\times_{\underline N^\gp/\underline N^\times}
\underline N^\sharp=\underline N
\end{equation}
(here it is needed that $\underline N$ is integral) such that the
diagram :
$$
\xymatrix{ P_{\!Y} \ar[r] \ar[d]_{\gamma_Y} & Q_Y \ar[r] \ar[d] &
P'_{\!Y} \ar[d]^{\alpha_Y} \\ \underline N \ar@{=}[r] &
\underline N \ar[r] & \underline N^\sharp }
$$
commutes. Then \eqref{eq_from-Q} determines a morphism of log
schemes $h_Y:(Y,\underline N)\to\Spec(\Z,Q)$, such that the
pair $(h_Y,\phi')$ is the unique morphism $h$ in $\cK_\intg$ with
the sought properties.

(iii): If both $P$ and $P'$ are fine, so is $Q$ (corollary
\ref{cor_fibres-are-fg}), and by construction, $\lambda^\gp$
is an isomorphism. Then the assertion follows from theorem
\ref{th_charact-smoothness}.
\end{proof}

\begin{proposition}\label{prop_basis-technique}
Let $\underline X:=((X,\underline M),F,\psi)$ be an object
of\/ $\cK_\intg$. Let also $\phi:F'\to F$ be an integral partial
subdivision, with $F$ locally fine and $F'$ integral. We have :
\begin{enumerate}
\item
$(X,\underline M)$ is a fine log scheme.
\item
If $F$ is saturated, $(X,\underline M)$ is a fs log scheme.
\item
$\underline X$ admits an inverse image over\/ $\phi$
(relative to the functor of \eqref{subsec_not_quite_fibration}).
\item
If $\phi$ is finite, then the cartesian morphism
$(f,\phi):\phi^*\underline X\to\underline X$, is quasi-compact.
\item
If $F'$ is locally fine, the morphism $(f,\phi)$ is quasi-separated
and \'etale.
\end{enumerate}
\end{proposition}
\begin{proof}(i): This is just a restatement of lemma
\ref{lem_soon-to-be-named}.

(ii): In light of (i), we only have to show that $\underline M_x$
is a saturated monoid, for every $x\in X$. Since by assumption
$\underline M_x^\sharp$ is saturated, the assertion follows from
lemma \ref{lem_exc-satura}(ii).

\begin{claim}\label{cl_reduce-t-aff-fan}
In order to show (iii)--(v), we may assume that :
\begin{enumerate}
\alphaenu
\item
$F=(\Spec\,P)^\sharp$ for a fine monoid $P$, and $X$ is an
affine scheme.
\item
The map of global sections
$P\to\Gamma(X,\underline M^\sharp)$ determined
by $\psi$, comes from a morphism of pre-log structures
$\beta:P_{\!X}\to\underline M$ which is a fine chart for
$\underline M$.
\end{enumerate}
\end{claim}
\begin{pfclaim} To begin with, suppose that
$((X',\underline M'),F',\psi')$ is the sought preimage of
$\underline X$; let $U\subset X$ be any open subset, and
$V\subset F$ an open subset such that $\psi(U)\subset V$.
Then it is easily seen that the object
$$
\phi^*\underline X\times_{\underline X}(U,V):=
((f^{-1}U,\underline M'_{|f^{-1}U}),\phi^{-1}V,\psi'_{|f^{-1}U})
$$
is a preimage of $((U,\underline M_{|U}),V,\psi_{|U})$ over the
restriction $\phi^{-1}V\to V$ of $\phi$. Now, suppose that we have
found an affine open covering $X=\bigcup_{i\in I}U_i$, and for every
$i\in I$ an affine open subset $V_i\subset F$ with $\psi(U_i)\subset
V_i$, such that the object $\underline U_i:=((U_i,\underline
M_{|U_i}),V_i,\psi_{|U_i})$ admits a preimage over the restriction
$\phi_i:\phi^{-1}V_i\to V_i$ of $\phi$; then the foregoing implies
that there are natural isomorphisms:
$$
\phi^*_i\underline U_i\times_{\underline U_i}(U_{ij},V_{ij})
\isom
\phi^*_j\underline U_j\times_{\underline U_j}(U_{ij},V_{ij})
$$
for every $i,j\in I$, where $U_{ij}:=U_i\cap U_j$ and
$V_{ij}:=V_i\cap V_j$. Thus, we may glue all these inverse images
along these isomorphisms, to obtain the sought inverse image of
$\underline X$.

Moreover, if $\phi$ is finite, the same will hold for the
restrictions $\phi_i$, and if each cartesian morphism
$\phi^*_i\underline U_i\to\underline U_i$ is quasi-compact, the same
will hold also for $(f,\phi)$ (\cite[Ch.I, \S6.1]{EGAI-new}).
Likewise, if each morphism $\phi^*_i\underline U_i\to\underline U_i$
is \'etale and quasi-separated, then the same will hold for
$(f,\phi)$ (\cite[Ch.I, Prop.6.1.11]{EGAI-new} and proposition
\ref{prop_sorite-smooth}(iii)).

Therefore, we may replace $X$ by any $U_i$, and $F$ by the
corresponding $V_i$, which reduces the proof of (iii)--(v)
to the case where condition (a) is fulfilled. Lastly, in light of
lemma \ref{lem_soon-to-be-named} we may suppose that the open subsets
$U_i$ are small enough, so that also condition (b) is fulfilled.
\end{pfclaim}

In view of claim \ref{cl_reduce-t-aff-fan}, we shall assume
henceforth that conditions (a) and (b) are fulfilled.
Let $\underline S:=(\Spec(\Z,P),T_P,\psi_P)$ be the object of
$\cK_\intg$ considered in example \ref{ex_fan-andlogscheme}(i) (with
$R:=\Z$); in this situation, $\beta$ determines a morphism of schemes
$f_\beta:X\to S$, and in view of lemma \ref{lem_localize-const} we
have $\underline X\simeq X\times_S\underline S$  (notation of
\eqref{subsec_base-change-in-cK}). Therefore, if we find an
inverse image $\underline S'$ for $\underline S$ over $\phi$, the object
$X\times_S\underline S'$ will provide the sought inverse image of 
$\underline X$. Thus, in order to show (iii)--(v), we may further reduce
to the case where $\underline X=\underline S$ (\cite[Ch.I, Prop.6.1.5(iii),
Prop.6.1.9(iii)]{EGAI-new} and proposition \ref{prop_sorite-smooth}(ii)).

\begin{claim}\label{cl_next-redux}
In order to prove (iii)--(v), we may assume that $F'$ is affine.
\end{claim}
\begin{pfclaim}
Indeed, say that $F'=\bigcup_{i\in I}V_i$ is an open covering,
and let $\phi_i:V_i\to T_P$ be the restriction of $\phi$; suppose
that we have found, for each $i\in I$, an inverse image
$\phi_i^*\underline S$ of $\underline S$ over $\phi_i$; again,
it follows easily that an inverse image for $\underline S$
over $\phi$ can be constructed by gluing the objects
$\phi_i^*\underline S$. This already implies that, in order
to show (iii), we may assume that $F'$ is affine.

Now, suppose that $\phi$ is finite, so we may find an open covering
as above, such that furthermore each $V_i$ is affine, and $I$ is a
finite set. Suppose that each of the corresponding morphisms
$\phi_i^*\underline S\to\underline S$ is quasi-compact; by the
foregoing, the schemes underlying the objects $\phi_i^*\underline S$
give a finite open covering of the scheme underlying
$\phi^*\underline S$, and then it is clear that (iv) holds.

Next, suppose furthermore that each morphism $\phi^*_i\underline
S\to\underline S$ is \'etale. It follows that the morphism
$\phi^*\underline S\to\underline S$ is also \'etale (proposition
\ref{prop_sorite-smooth}(ii)); then, since $S$ is noetherian, the
scheme underlying $\phi^*\underline S$ is locally noetherian
(\cite[Ch.I, Prop.6.2.2]{EGAI-new}), and therefore the morphism
$\phi^*\underline S\to\underline S$ is quasi-separated (\cite[Ch.I,
Cor.6.1.13]{EGAI-new}).
\end{pfclaim}

In view of claim \ref{cl_next-redux}, we shall assume that
$F'=(\Spec\,P')^\sharp$ is affine as well, for a sharp and
integral $P'$, and that $\phi$ is given by a morphism
$\lambda:P\to P'$ inducing a surjection on the associated
groups (details left to the reader). Then assertions (iii)
and (iv) are now straightforward consequences of lemma
\ref{lem_redux-more}(i,ii), and (v) follows from lemma
\ref{lem_redux-more}(iii), after one remarks that, when $F'$ is
fine, one may choose for $P'$ a fine monoid.
\end{proof}

\begin{remark}\label{rem_basic-technique}
Keep the assumptions of proposition \ref{prop_basis-technique},
and let $t\in F_0$ a point of $F$ of height zero
(see \eqref{subsec_height-in-T}). Notice that the inclusion map
$j_t:\{t\}\to F$ is an open immersion, hence the fibre
$X_t:=\psi^{-1}(t)\subset X$ is open; indeed, it is clear from
the definitions, that $\psi^{-1}(F_0)$ is precisely the trivial
locus of $\underline X$. Moreover, let $t'\in\phi^{-1}(t)$ since
the group homomorphism $\cO^\gp_{\!F,t}\to\cO^\gp_{\!F',t'}$
is surjective, we see that $t'$ is of height zero in $F'$, and
$\phi$ restricts to an isomorphism of fans
$(\{t'\},\cO_{\!F,t'})\isom(\{t\},\cO_{\!F,t})$.
Set $\underline X{}_t:=j_t^*\underline X$ (whose underlying scheme
is $X_t$), and define likewise
$\phi^*\underline X{}_{t'}:=(\phi\circ j_{t'})^*\underline X$
(whose underlying scheme is an open subset of the trivial locus
of $\phi^*\underline X$).  We deduce that :
\begin{itemize}
\item
The trivial locus of $\phi^*\underline X$ is the preimage of
$(X,\underline M)_\tr$.
\item
For every $t'\in F'_0$, the restriction of $(f,\phi)$ :
$\phi^*\underline X{}_{t'}\to\underline X{}_t$ is an isomorphism.
\end{itemize}
\end{remark}

\begin{corollary}\label{cor_case-of-blow-sat}
In the situation of proposition {\em\ref{prop_basis-technique}},
let $((X',\underline M'),F',\psi'):=\phi^*\underline X$.
The following holds :
\begin{enumerate}
\item
If $F'=F^\sat$, and $\phi:F^\sat\to F$ is the counit of adjunction,
then $(X',\underline M')=(X,\underline M)^\mathrm{fs}$, and $f$ is
the counit of adjunction.
\item
If $\phi$ is the blow up of a coherent ideal $\cI$ of $\cO_F$, then
$f$ is the blow up of $\cI\!\underline M$, the unique ideal of
$\underline M$ whose image in $\underline M^\sharp$ equals
$\psi^{-1}\!\cI$.
\item
Suppose moreover, that $(X,\underline M)$ is regular, $F'$ is
saturated , and $\underline X=\cK(X,\underline M)$ (notation
of remark {\em\ref{rem_discuss-fan-logreg}(iii)}). Then
$(X',\underline M')$ is regular, and $\phi^*\underline X$ is
isomorphic to $\cK(X',\underline M')$.
\end{enumerate}
\end{corollary}
\begin{proof} To start with, we remark that the assertions are
local on $X$. Indeed, this is clear for (i), and for (iii) it
follows easily from remark \ref{rem_discuss-fan-logreg}(i,,ii);
concerning (ii), suppose that $X=\bigcup_{i\in I}U_i$ is an open
covering, such that $\phi^*(U_i,\underline M_{|U_i})$ is the blow
up of the ideal $\cI\!\underline M_{|U_i}$, for every $i\in I$.
For every $i,j\in I$ set $U_{ij}:=U_i\cap U_j$; by the universal
property of the blow up, there are unique isomorphisms of
$(U_{ij},\underline M_{|U_{ij}})$-schemes :
$$
U_{ij}\times_{U_i}\phi^*(U_i,\underline M_{|U_i})\isom
U_{ij}\times_{U_j}\phi^*(U_j,\underline M_{|U_j}).
$$
Then both $\phi^*(X,\underline M)$ and the blow up of
$\cI\!\underline M$ are necessarily obtained by gluing along
these isomorphisms, so they are isomorphic.

Thus, we may assume that $F=(\Spec\,P)^\sharp$ for some
fine monoid $P$, and $(X,\underline M)=X\times_S(S,P^{\log}_S)$
for some morphism of schemes $X\to S$, where as usual
$(S,P^{\log}_S,F,\psi_P)$ is defined as in example
\ref{ex_fan-andlogscheme}(i). In this case, (i) follows
from remarks \ref{rem_quasi-fine}(iii),
\ref{rem_product-of-aff-fans}(ii) and lemma
\ref{lem_redux-more}(ii).

Likewise, remark \ref{rem_quite}(ii) allows to reduce (ii)
to the case where $\underline X=(S,P^{\log}_S,F,\psi_P)$,
and $\cI=I^\sim$ for some ideal $I\subset P$; in which case
we conclude by inspecting the explicit description in example
\ref{ex_explicit-blowup}.

(iii): From proposition \ref{prop_basis-technique}(ii,v)
and theorem \ref{th_smooth-preserve-reg} we already see
that $(X',\underline M')$ is regular. Next, we are easily
reduced to the case where $X$ is affine, say $X=\Spec\,A$, and
$F'=(\Spec\,Q)^\sharp$, for some saturated monoid $Q$, and
by lemma \ref{lem_redux-more}, we may assume that $\phi$ is
induced by a morphism of monoids $\lambda:P\to Q$ such that
$\lambda^\gp$ is an isomorphism, and
$(X',\underline M')=$ $\phi^*\underline X=X\times_S\underline S'$,
where $\underline S':=(\Spec(\Z,Q),F',\psi_Q)$ is defined as
in example \ref{ex_fan-andlogscheme}(ii). Let $\fq\subset Q$
be any prime ideal, and set $\fp:=\fq\cap P$; in light of
remark \ref{rem_discuss-fan-logreg}(i), it then suffices to
show that $Q/\fq\otimes_PA=Q/\fq\otimes_{P/\fp}A/\fp A$ is an
integral domain. To this aim, let us remark :

\begin{claim}\label{cl_free-poly}
Let $\lambda:P\to Q$ be a morphism of fine monoids, such that
$\lambda^\gp$ is an isomorphism, $F\subset Q$ any face, and
denote $\lambda_F:F\cap P\to P$ the inclusion map. Then :
\begin{enumerate}
\item
The natural map $\Coker\,\lambda^\gp_F\to P^\gp/(F\cap P)^\gp$
is injective.
\item
If moreover, $P$ is saturated, then $\Coker\,\lambda^\gp_F$ is
a free abelian group of finite rank.
\end{enumerate} 
\end{claim}
\begin{pfclaim} The map of (i) is obtained via the snake
lemma, applied to the ladder of abelian groups :
$$
\xymatrix{
0 \ar[r] & (F\cap P)^\gp \ar[r] \ar[d]_{\lambda_F^\gp} &
P^\gp \ar[r] \ar[d]^{\lambda^\gp} &
P^\gp/(F\cap P)^\gp \ar[r] \ar[d] & 0 \\
0 \ar[r] & F^\gp \ar[r] & Q^\gp \ar[r] & Q^\gp/F^\gp \ar[r] & 0.
}$$
taking into account that both $\Ker\,\lambda^\gp$ and
$\Coker\,\lambda^\gp$ vanish. Then (i) is obvious, and
(ii) comes down to checking that $P^\gp/(F\cap P)^\gp$
is torsion-free, in case $P$ is fine and saturated.
But since $F\cap P$ is a face of $P$, the latter assertion
is an easy consequence of proposition \ref{prop_Gordie}.
\end{pfclaim}

Now, set $F:=Q\setminus\fq$; we come down to checking
that $B:=F\otimes_{F\cap P}A/\fp A$ is an integral domain,
and remark \ref{rem_discuss-fan-logreg}(i) tells us that
$A/\fp A$ is a domain. However, $A/\fp A$ is
$(F\cap P)$-flat (proposition \ref{prop_second-crit}(ii)),
hence the natural map
$$
B\to C:=F^\gp\otimes_{(F\cap P)^\gp}\Frac(A/\fp A)
$$
is injective; on the other hand, claim \ref{cl_free-poly}(ii)
implies that $C=A[\Coker\,\lambda^\gp_F]$ is a free
(polynomial) $A$-algebra, whence the contention.
\end{proof}

\begin{example}\label{ex_fan-andlogsch-2}
In the situation of example \ref{ex_uniqueness-of-mu_n}, take
$T:=\Spec\,P$ for a fine monoid $P$, and define $\underline S$
as in example \ref{ex_fan-andlogscheme}(i). The $k$-Frobenius
map $\bek_P$ (example \ref{ex_multiply-by-n-in-fan}(i)) induces
an endomorphism $\bek_{\underline S}:=(\Spec(R,\bek_P),\bek_T)$
of $\underline S$ in $\cK$. By proposition
\ref{prop_basis-technique}(iii) and example \ref{ex_uniqueness-of-mu_n},
there exists a unique morphism $\underline g$ fitting into a
commutative diagram of $\cK_\intg$ :
\set\begin{equation}\label{eq_yad-in-cK}
{\diagram 
\phi^*\underline S \ar[r] \ar[d]_{\underline g} &
\underline S \ar[d]^{\bek_{\underline S}} \\
\phi^*\underline S \ar[r] & \underline S
\enddiagram}
\end{equation}
(where the horizontal arrows are the cartesian morphisms). Say that
$\phi^*\underline S=((Y,\underline N),F,\psi)$; then we have
$\underline g=(g,\bek_F)$ for a unique endomorphism $g:=(g,\log g)$
of the log scheme $(Y,\underline N)$. Let $U:=\Spec\,P'\subset F$ be
any open affine subset ; since $\bek_F$ is the identity on the
underlying topological spaces, $g$ restricts to an endomorphism
$g_{|\psi^{-1}U}$ of $\psi^{-1}U\times_Y(Y,\underline N)$. In view
of lemma \ref{lem_redux-more}, the latter log scheme is of the form
$\underline S'$ as in example \ref{ex_fan-andlogscheme}(ii), with
$Q:=P^\gp\times_{P^{\prime\gp}}P'$. Then $g_{|\psi^{-1}U}$
is induced by an endomorphism $\nu$ of $Q$, fitting into a commutative
diagram :
$$
\xymatrix{ P \ar[r] \ar[d]_{\bek_P} & Q \ar[d]^\nu \\
                  P \ar[r] & Q 
}$$
whose horizontal arrows are the natural injections. Since $Q\subset P^\gp$,
it is clear that $\nu=\bek_Q$. Especially, $g:Y\to Y$ is a finite morphism
of schemes. Furthermore, for every point $y\in Y$, we have a commutative
diagram of monoids :
\set\begin{equation}\label{eq_finnegan}
{\diagram
\underline N{}_{g(y)} \ar[d]_{\log g_y} \ar[r] & 
\underline N{}^\sharp_{g(y)} \ar[d]_{\log g^\sharp_y} &
\cO_{F,\psi(y)} \ar[l]_-\sim
\ar[d]^{\bek_{\psi(y)}:=(\log\bek_F)_{\psi(y)}} \\
\underline N{}_y \ar[r] & \underline N{}^\sharp_y &
\cO_{F,\psi(y)} \ar[l]_-\sim.
\enddiagram}
\end{equation}
\end{example}

\sset\subsubsection{}\label{subsec_canon-obj-valuat}
Let $(K,|\cdot|)$ be a valued field, $\Gamma$ and $K^+$
respectively the value group and the valuation ring of $|\cdot|$;
denote by $1\in\Gamma$ the neutral element, and by
$\Gamma_{\!+}\subset\Gamma$ the submonoid consisting of all
elements $\leq 1$. Set $S:=\Spec\,K^+$; we consider
the log scheme $(S,\cO_{\!S}^*)$, where
$\cO_{\!S}^*\subset\cO_{\!S}$ is the subsheaf such that
$\cO_{\!S}^*(U):=\cO_{\!S}(U)\!\setminus\!\{0\}$ for every open
subset $U\subset S$. To this log scheme we associate the object
of $\cK_\intg$ :
$$
\sff(K,|\cdot|):=((S,\cO^*_{\!S}),\Spec\,\Gamma_{\!+},\psi_\Gamma)
$$
where $\psi_\Gamma$ is the morphism of monoidal spaces arising from
the isomorphism
$$
\Gamma_{\!+}\isom(K^+\!\setminus\!\{0\})/(K^+)^\times
$$
deduced from the valuation $|\cdot|$. It is well known that
$\psi_\Gamma$ is a homeomorphism (see {\em e.g.}
\cite[\S6.1.26]{Ga-Ra}). We shall use objects of this kind to state
some separation and properness criteria for log schemes whose
existence is established via proposition \ref{prop_basis-technique}.
To this aim, we need to digress a little, to prove the following
auxiliary results, which are refinements of the standard valuative
criteria for morphisms of schemes.

\sset\subsubsection{}\label{subsec_max-points-in-img}
Let $f:X\to Y$ be a morphism of schemes, $R$ an integral ring, and
$K$ the field of fractions of $R$.
Let us denote by $X(R)_{\max}\subset X(R)$ the set of morphisms
$\Spec\,R\to X$ which map $\Spec\,K$ to a maximal point of $X$.
There follows a commutative diagram of sets :
\set\begin{equation}\label{eq_max-sections}
{\diagram
X(R)_{\max} \ar[r] \ar[d] & X(K)_{\max} \ar[d] \\
Y(R) \ar[r] & Y(K).
\enddiagram}
\end{equation}

\begin{proposition}\label{prop_separate-gen}
Let $f:X\to Y$ be a morphism of schemes. The following conditions
are equivalent :
\begin{enumerate}
\alphaenu
\item
$f$ is separated.
\item
$f$ is quasi-separated, and for every valued field $(K,|\cdot|)$,
the map of sets :
\set\begin{equation}\label{eq_map-of-max}
X(K^+)_{\max}\to X(K)_{\max}\times_{Y(K)}X(K^+)_{\max}
\end{equation}
deduced from diagram \eqref{eq_max-sections} (with $R:=K^+$, the
valuation ring of\/ $|\cdot|$), is injective.
\end{enumerate}
\end{proposition}
\begin{proof} (a) $\Rightarrow$ (b) by the valuative criterion
of separation (\cite[Ch.I, Prop.5.5.4]{EGAI-new}). Conversely, we
shall show that if (b) holds, then the assumptions of the criterion
of {\em loc.cit.} are fulfilled. Indeed, let $(L,|\cdot|_L)$ be any
valued field, with valuation ring $L^+$, and suppose we have two
morphisms $\sigma_1,\sigma_2:\Spec\,L^+\to X$ whose restrictions
to $\Spec\,L$ agree, and such that $f\circ\sigma_1=f\circ\sigma_2$.

Let $s,\eta\in\Spec\,L^+$ be respectively the closed point and the
generic point, set $x:=\sigma_1(\eta)=\sigma_2(\eta)\in X$, and
denote by $\phi:\cO_{\!X,x}\to L$ the ring homomorphism corresponding
to the restriction of $\sigma_1$ (and $\sigma_2$); then $x$ admits
two specializations $x_i:=\sigma_i(s)\in X$ (for $i=1,2$) such that
$\phi$ sends the image $A_i\subset\cO_{\!X,x}$ of the specialization
map $\cO_{\!X,x_i}\to\cO_{\!X,x}$, into $L^+$, and the maximal ideal
of $A_i$ into the maximal ideal $\fm_L$ of  $L^+$, for both $i=1,2$.
Moreover, $y:=f(x_1)=f(x_2)$.

Denote by $B\subset\cO_{\!X,x}$ the smallest subring containing
$A_1$ and $A_2$, and set $\fp:=B\cap\phi^{-1}\fm_L$. Then
$\phi(B_\fp)\subset L^+$ as well, and $B_\fp$ dominates both $A_1$
and $A_2$. Now, let $t\in X$ be a maximal point which specializes
to $x$; by \cite[Ch.I, Prop.5.5.2]{EGAI-new} we may find a valued
field $(E,|\cdot|_E)$ with a local ring homomorphism $\cO_{\!X,t}\to E$,
such that the valuation ring $E^+$ of $|\cdot|_E$ dominates
the image of the specialization map $\cO_{\!X,x}\to\cO_{\!X,t}$.
Let $\kappa(E)$ be the residue field of $E^+$, and
$\bar B_\fp\subset\kappa(E)$ the image of $B_\fp$. By the same
token, we may find a valuation ring $V\subset\kappa(E)$ with
fraction field $\kappa(E)$, and dominating $\bar B_\fp$; then
it is easily seen that the preimage $K^+\subset E^+$ of $V$ is
a valuation ring with field of fractions $K=E$, which dominates the
image of $B_\fp$ in $\cO_{\!X,t}$ (\cite[Th.10.1(iv)]{Mat}).
Hence, $K^+$ dominates the images of $\cO_{\!X,x_i}$, for
both $i=1,2$, and therefore, also the image of $\cO_{\!Y,y}$;
in other words, in this way we obtain two elements
in $X(K^+)_{\max}$ whose images agree in $Y(K^+)$, and whose
restrictions $\Spec\,K\to X$ coincide. By assumption, these
two $K^+$-points must then coincide, especially $x_1=x_2$, and
therefore $\sigma_1=\sigma_2$, as required.
\end{proof}

\begin{proposition}\label{prop_univ-close-gen}
Let $f:X\to Y$ be a quasi-compact morphism of schemes.
The following conditions are equivalent :
\begin{enumerate}
\alphaenu
\item
$f$ is universally closed.
\item
For every valued field $(K,|\cdot|)$, the corresponding map
\eqref{eq_map-of-max} is surjective.
\end{enumerate}
\end{proposition}
\begin{proof} (a) $\Rightarrow$ (b) by the valuative criterion
of \cite[Ch.I, Prop.5.5.8]{EGAI-new}. Conversely, we will show
that (b) implies that the conditions of {\em loc.cit.} are
fulfilled. Hence, let $(L,|\cdot|_L)$ be any valued field,
with valuation ring $L^+$, and suppose that we have a morphism
$\sigma:\Spec\,L\to X$, whose composition with $f$ extends to
a morphism $\Spec\,L^+\to Y$; we have to show that $\sigma$
extends to a morphism $\Spec\,L^+\to X$. Denote by
$\eta,s\in\Spec\,L^+$ respectively the generic and the closed
point, let $x\in X$ be the image of $\eta$, and $y\in Y$ the
image of $s$. Then $\sigma$ corresponds to a ring homomorphism
$\sigma^\natural:\cO_{\!X,x}\to L$, and $L^+$ dominates the
image of the map $\cO_{Y,y}\to\cO_{\!X,x}$ determined by $f$.
Moreover, $\sigma^\natural$ factors through the residue field
$\kappa(x)$ of $\cO_{\!X,x}$. Then $\kappa(x)^+:=\kappa(x)\cap L^+$
is a valuation ring with fraction field $\kappa(x)$, and we are
reduced to showing that there exists a specialization $x'\in X$
of $x$, such that $f(x')=y$, and such that $\kappa(x)^+$ dominates
the image of the specialization map $\cO_{\!X,x'}\to\cO_{\!X,x}$.

Let now $t\in X$ be a maximal point which specializes to $x$;
by \cite[Ch.I, Prop.5.5.8]{EGAI-new} we may find a valued field
$(E,|\cdot|_E)$ with a local ring homomorphism $\cO_{\!X,t}\to E$,
whose valuation ring $E^+$ dominates the image of the specialization
map $\cO_{\!X,x}\to\cO_{\!X,t}$. Let $\kappa(E)$ be the residue
field of $E^+$; the induced map $\cO_{\!X,x}\to\kappa(E)$ factors
through $\kappa(x)$, and by \cite[Th.10.2]{Mat} we may find a
valuation ring $V\subset\kappa(E)$ with field of fractions
$\kappa(E)$, which dominates $\kappa(x)^+$. The preimage
$K^+\subset E^+$ of $V$ is a valuation ring with field of
fractions $K=E$ (\cite[Th.10.1]{Mat}). By construction, $K^+$
dominates the image of $\cO_{\!X,y}$, in which case assumption
(b) says that there exists a specialization $x'$ of $x$ such that
$f(x')=y$, and such that $K^+$ dominates the image of the
specialization map $\cO_{\!X,x'}\to\cO_{\!X,t}$. A simple
inspection then shows that $\kappa(x)^+$ dominates the image
of $\cO_{\!X,x'}$ in $\cO_{\!X,x}$, as required.
\end{proof}

\begin{corollary}\label{cor_proper-crit-fanslog}
Let $f:X\to Y$ be a morphism of schemes which is quasi-separated
and of finite type. Then the following conditions
are equivalent :
\begin{enumerate}
\alphaenu
\item
$f$ is proper.
\item
For every valued field $(K,|\cdot|)$, the corresponding diagram
\eqref{eq_max-sections} (with $R:=K^+$ the valuation ring of\/
$|\cdot|$) is cartesian.
\end{enumerate}
\end{corollary}
\begin{proof} It is immediate from propositions \ref{prop_separate-gen}
and \ref{prop_univ-close-gen}.
\end{proof}

\sset\subsubsection{}\label{subsec_almost-home}
We are now ready to return to log schemes. Resume the situation
of \eqref{subsec_canon-obj-valuat}, let
$\underline X:=((X,\underline M),F,\psi)$ be any object of $\cK$,
and denote by $\alpha:\underline M\to\cO_{\!X}$ the structure map of
$\underline M$.

Suppose we are given a morphism $\sigma:S\to X$ of schemes, and we
ask whether there exists a morphism of log structures
$\beta:\sigma^*\underline M\to\cO^*_{\!S}$, such that the pair
$(\sigma,\beta)$ is a morphism of log schemes
$(S,\cO^*_{\!S})\to(X,\underline M)$. By definition, this holds if
and only if the composition
$$
\bar\beta:\sigma^*\underline M\xrightarrow{ \sigma^*\alpha }
\sigma^*\cO_{\!X}\xrightarrow{ \sigma^\natural }\cO_{\!S}
$$
factors through $\cO^*_{\!S}$. Moreover, in this case the
factorization is unique, so that $\sigma$ determines $\beta$
uniquely. Let $\eta\in S$ be the generic point; we claim that the
stated condition is fulfilled, if and only if $t:=\sigma(\eta)$
lies in $(X,\underline M)_\tr$ (see definition
\ref{def_trivial-locus}(i)).
Indeed, if $t$ lies in the trivial locus,
$\underline M_t=\cO_{\!X,t}^\times$, so certainly the image
of $\bar\beta_t:\underline M_t\to\cO_{\!S,\eta}$ lies in
$\cO_{\!S,\eta}^*=K^\times$. Now, if $s\in S$ is any other
point, the composition
$$
\underline M_{\sigma(s)}\xrightarrow{ \bar\beta_{\sigma(s)} }
\cO_{\!S,s}\to\cO_{\!S,\eta}=K
$$
factors through $\underline M_t$, hence its image lies in
$K^\times\cap\cO_{\!S,s}=\cO^*_{\!S,s}$, whence the contention.
Conversely, if $\bar\beta$ factors through $\cO^*_{\!S}$,
it follows especially that the image of the stalk of the
structure map $\underline M_t\to\cO_{\!X,t}$ lies in the
preimage of $\cO^*_{\!S,\eta}=\cO^\times_{\!S,\eta}$, and
the latter is just $\cO_{\!X,t}^\times$, so $t$ lies in the
trivial locus.

Next, let $f:(S,\cO_{\!S}^*)\to(X,\underline M)$ be any morphism of
log schemes. We claim that there exists a unique morphism of fans
$\phi:=(\phi,\log\phi):\Spec\,\Gamma_+\to F$, such that the pair
$(f,\phi)$ is a morphism $\sff(K,|\cdot|)\to\underline X$ in $\cK$.
Indeed, since $\psi_\Gamma$ is a homeomorphism, there exists a
unique continuous map $\phi$ on the underlying topological spaces,
such that $\phi\circ\psi_\Gamma=\psi\circ f$, and for the same
reason, the map $\log f^\sharp:f^*\underline
M/\cO^\times_{\!S}\to\cO^*_{\!S}/\cO_{\!S}^\times$ is of the form
$\psi_\Gamma^*(\log\phi)$ for a unique morphism of sheaves
$\log\phi$ as sought. Then $\log\phi$ will be a local morphism,
since the same holds for $\log f$ (see
\eqref{subsec_resolve-logscheme}).

Summing up, we have shown that the natural map
$$
\Hom_\cK(\sff(K,|\cdot|),\underline X)\to X(K^+) \quad : \quad
(\sigma,\phi)\mapsto\sigma
$$
is injective, and its image is the set of all the morphisms $S\to X$
which map $\eta$ into $(X,\underline M)_\tr$.

\begin{proposition}\label{prop_proper-crit-in-C}
Let $\underline X:=((X,\underline M),F,\psi)$ be an object
of\/ $\cK_\intg$, and $\phi:F'\to F$ an integral partial
subdivision, with $F$ and $F'$ locally fine. We have :
\begin{enumerate}
\item
If the induced map $F'(\N)\to F(\N)$ is injective, the cartesian
morphism $\phi^*\underline X\to\underline X$ is separated.
\item
If $\phi$ is a proper subdivision, the cartesian morphism
$\phi^*\underline X\to\underline X$ is proper, and induces
an isomorphism of schemes
$(\phi^*\underline X)_\tr\isom(X,\underline M)_\tr$.
\end{enumerate}
\end{proposition}
\begin{proof} The assertions are local on $X$ (cp. the proof of
claim \ref{cl_reduce-t-aff-fan}), so we may assume -- by lemma
\ref{lem_soon-to-be-named} -- that $\underline M$ admits
a chart $P_{\!X}\to\underline M$, and $F=(\Spec\,P)^\sharp$.
In this case, let $S:=\Spec\,\Z[P]$, and denote by $\underline S$
the object of $\cK_{\intg}$ attached to $P$, as in example
\ref{ex_fan-andlogscheme}(i) (with $R:=\Z$); then
$(X,\underline M)$ is isomorphic to $X\times_S(S,P^{\log}_S)$,
and if $f_S:\phi^*\underline S\to\underline S$ is the cartesian
morphism over $\phi$, then the cartesian morphism
$f:\phi^*\underline X\to\underline X$ is given by the pair
$(\one_X\times_Sf_S,\phi)$ (see the proof of proposition
\ref{prop_basis-technique}(iii)). Thus, we may replace
$\underline X$ by $\underline S$
(\cite[Ch.I, Prop.5.3.1(iv)]{EGAI-new}), in which case lemma
\ref{lem_redux-more} shows that the log scheme
$(X',\underline M')$ underlying $\phi^*\underline S$ admits an
open covering consisting of affine log schemes of the form
$\Spec(\Z,Q)$, for a fine monoid $Q$. Notice that, for such
$Q$ we have :
$$
\Spec(\Z,Q)_\tr=\Z[Q^\gp]
$$
which is a dense open subset of $\Spec\,\Z[Q]$.

(i): According to proposition \ref{prop_basis-technique}(v), the
morphism $f$ is quasi-separated, so we may apply the criterion of
proposition \ref{prop_separate-gen}. Indeed, let $(K,|\cdot|)$ be
any valued field, and suppose that $\sigma_i\in X'(K^+)_{\max}$
(for $i=1,2$) are two sections, whose images in
$X'(K)_{\max}\times_{S(K)}S(K^+)$ coincide; we have to show that
$\sigma_1=\sigma_2$. However, we have just seen that the maximal
points of $X'$ lie in $(X',\underline M')_\tr$, so the discussion
of \eqref{subsec_almost-home} shows that both $\sigma_i$ extend
uniquely to morphisms
$\sigma'_i:\sff(K,|\cdot|)\to\phi^*\underline S$ in $\cK$, and it
suffices to show that $\sigma'_1=\sigma'_2$. By definition, the
datum of $\sigma'_i$ is equivalent to the datum of a morphism
$\sigma''_i:\sff(K,|\cdot|)\to\underline S$, and a morphism of
fans $\phi'_i:\Spec\,\Gamma_{\!+}\to F'$. Again by
\eqref{subsec_almost-home}, the morphisms $\sigma''_1$ and
$\sigma''_2$ agree if and only if they induce the same morphisms
of schemes; the latter holds by assumption, since $\sigma_1$ and
$\sigma_2$ yield the same element of $S(K^+)$. On the other hand,
in view of (b) and proposition \ref{prop_val-crit-fans},
the elements $\phi_1,\phi_2\in F'(\Gamma_{\!+})$ coincide if
and only if their images in $F(\Gamma_{\!+})$ coincide; but
again, this last condition holds since the images of both
$\sigma_i$ agree in $S(K^+)$.

(ii): In view of (i) and proposition
\ref{prop_basis-technique}(iv),(v) we know already that $f$ is
separated and of finite type, so we may apply the criterion of
corollary \ref{cor_proper-crit-fanslog}. Hence, let $\sigma\in
S(K^+)$ be a section, and $x\in X'(K)_{\max}$ a $K$-rational point
such that $f(x)$ is the image of $\sigma$ in $S(K)$; in view of (i),
it suffices to show that $\sigma$ lifts to a section
$\tilde\sigma\in X'(K^+)_{\max}$, whose image in $X'(K)$ is $x$.
Since $x\in(X',\underline M')_\tr$, remark \ref{rem_basic-technique}
implies that $f(x)\in(S,P^{\log}_{\!S})_\tr$, and then the
discussion in \eqref{subsec_almost-home} says that $\sigma$
underlies a unique morphism
$(\sigma',\beta):\sff(K,|\cdot|)\to\underline S$. By proposition
\ref{prop_order-only}, the element $\beta\in F(\Gamma_{\!+})$ lifts
to an element $\beta'\in F'(\Gamma_{\!+})$, and finally the pair
$((\sigma,\beta),\beta')$ determines a unique morphism
$\sff(K,|\cdot|)\to\phi^*\underline S$, whose underlying morphism of
schemes is the sought $\tilde\sigma$.
Lastly, notice that the map $F'(\{1\})\to F(\{1\})$
induced by $\phi$, is bijective by propositions
\ref{prop_val-crit-fans} and \ref{prop_order-only} (where $\{1\}$
is the monoid with one element), which means that $\phi$
restricts to a bijection on the points of height zero; then
remark \ref{rem_basic-technique} implies the second assertion
of (ii).
\end{proof}

\begin{theorem}\label{th_resolution-Zar}
Let $(X,\underline M)$ be a regular log scheme. Then there
exists a smooth morphism of log schemes
$f:(X',\underline M')\to(X',\underline M)$, whose underlying
morphism of schemes is proper and birational, and such that
$X'$ is a regular scheme. More precisely, $f$ restricts to an
isomorphism of schemes $f^{-1}X_\reg\to X_\reg$ on the preimage
of the open locus of regular points of $X$.
\end{theorem}
\begin{proof} We use the object
$\underline X:=((X,\underline M),F(X,\underline M),\pi_X)$ attached
to $(X,\underline M)$ as in remark \ref{rem_discuss-fan-logreg}(ii).
Indeed, it is clear that $F(X,\underline M)$ is locally fine and
saturated, hence theorem \ref{th_Saint-Donat} yields an integral,
proper, simplicial subdivision $\phi:F'\to F(X,\underline M)$
which restricts to an isomorphism
$\phi^{-1}F(X,\underline M)_\mathrm{sim}\isom
F(X,\underline M)_\mathrm{sim}$. Take
$(f,\phi):\phi^*\underline X\to\underline X$ to be the cartesian
morphism over $\phi$, and denote by $(X',\underline M')$ the log
scheme underlying $\phi^*\underline X$; it follows already from
proposition \ref{prop_proper-crit-in-C}(ii) that $f$ is proper on
the underlying schemes. Next, corollary \ref{cor_more-precisely}
shows that $X_\mathrm{reg}$ is
$\pi^{-1}F(X,\underline M)_\mathrm{sim}$, so $f$ restricts
to an isomorphism $f^{-1}X_\mathrm{reg}\isom X_\mathrm{reg}$.
Furthermore, $f$ is \'etale, by proposition
\ref{prop_basis-technique}(v), hence the log scheme
$(X',\underline M')$ is regular
(theorem \ref{th_smooth-preserve-reg}). Finally, again by corollary
\ref{rem_more-precisely} we see that $X'$ is regular.
\end{proof}

\sset\subsubsection{}
Let now $(Y,\underline N)$ be a regular log scheme, such that
$F(Y,\underline N)$ is affine (notation of remark
\ref{rem_discuss-fan-logreg}(ii)), say isomorphic to
$(\Spec\,P)^\sharp$, for some fine, sharp and saturated monoid
$P$. Let $I\subset P$ be an ideal generated by two elements
$a,b\in P$, and denote by
$f:(Y',\underline N')\to(Y,\underline N)$ the saturated blow
up of the ideal $I\underline N$ of $\underline N$ (see
\eqref{subsec_sat-blowup}). Set $U':=(Y',\underline N')_\tr$,
$U:=(Y,\underline N)_\tr$, and denote $j:U\to Y$, $j':U'\to Y'$
the open immersions. In this situation we have :

\begin{lemma}\label{lem_conclude-indu}
{\em(i)}\ \ $H^1(Y',\underline N'{}^{\sharp\gp})=0$.
\begin{enumerate}
\addenu
\item
Suppose morever, that $R^1j'_*\cO^\times_{\!U'}=0$. Then
$R^1j_*\cO^\times_{\!U}=0$.
\end{enumerate}
\end{lemma}
\begin{proof}(i): Since
$\underline N'{}^{\sharp\gp}=\pi_X^*\cO^\gp_{\!F(Y',\underline N')}$,
claim \ref{cl_open-irredu}(ii) and remark
\ref{rem_discuss-fan-logreg}(iii) reduce to showing
\set\begin{equation}\label{eq_vanish-this}
H^1(F(Y',\underline N'),\cO^\gp_{\!F(Y',\underline N')})=0.
\end{equation}
However, by corollary \ref{cor_case-of-blow-sat}(iii), the fan
$F(Y',\underline N')$ is the saturated blow up of the ideal
$I\cO_{\!F(Y,\underline N)}$ of $\cO_{\!F(Y,\underline N)}$,
hence it admits
the affine covering
$$
F(Y',\underline N')=
(\Spec\,P[a^{-1}b]^\sat)^\sharp\cup(\Spec\,P[b^{-1}a]^\sat)^\sharp.
$$
Notice now that every affine fan is a local topological space,
hence the left hand-side of \eqref{eq_vanish-this} is computed
by the \v{C}ech cohomology of $\cO^\gp_{\!F(Y',\underline N')}$
relative to this covering (theorem \ref{th_Cartan}(ii)).
However, the intersection of the two open subsets is
$(\Spec\,P[a^{-1}b,b^{-1}a])^\sharp$, and clearly the
restriction map
$$
H^0(\Spec\,P[a^{-1}b]^\sat,\cO^\gp_{F(Y',\underline N')})\to
H^0(\Spec\,P[a^{-1}b,b^{-1}a],\cO^\gp_{F(Y',\underline N')})
$$
is surjective. The assertion is an immediate consequence.

(ii): Let $y\in Y$ be any point; according to \eqref{eq_quick}, the
morphism
$$
f\times_YY(y):(Y',\underline N')\times_YY(y)\to(Y(y),\underline N(y))
$$
is the saturated blow up of the ideal $I\underline N(y)$; on
the other hand, let $U(y):=U\times_YY(y)$, $U'(y):=U'\times_YY(y)$
and denote by $j_y:U(y)\to Y(y)$ and $j'_y:U'(y)\to Y'\times_YY(y)$
the open immersions; in light of proposition
\ref{prop_dir-im-and-colim}(ii), it suffices to show that
$R^1j_{y*}\cO_{\!U(y)}=0$, and the assumption implies that
$R^1j'_{y*}\cO_{\!U'(y)}=0$. Summing up, we may replace
$Y$ be $Y(y)$, and assume from start that $Y$ is local, and
$y$ is its closed point. From the assumption we get :
$$
H^1(Y',j'_*\cO^\times_{\!U'})=H^1(U',\cO^\times_{Y'})=\Pic\,U'.
$$
On the other hand, recall that
$\underline N'{}^{\sharp\gp}=j'_*\cO^\times_{\!U'}/\cO^\times_{Y'}$
(proposition \ref{prop_log-struct-is-fixed}); combining
with (i), we deduce that the natural map
$$
\Pic\,Y'\to\Pic\,U'
$$
is surjective. Set $Y'_0:=f^{-1}(y)\subset Y'$, endow $Y'_0$
with its reduced subscheme structure, and let $i:Y'_0\to Y'$
be the closed immersion. If $I$ is an invertible ideal of $P$,
then $f$ is an isomorphism, in which case the assertion is
obvious. We may then assume that $I$ is not invertible,
in which case claim \ref{cl_more-and-more}(ii) says that
there exists a morphism of $(Y,\underline N)$-schemes
$h:(Y',\underline N')\to\P^1_{(Y,\underline N)}$
inducing an isomorphism of $\kappa(y)$-schemes
\set\begin{equation}\label{eq_this-is-h}
(h\times_Y\Spec\,\kappa(y))_\red:Y'_0\isom\P^1_{\kappa(y)}.
\end{equation}
Let us remark :

\begin{claim}\label{cl_from-thesis}
Let $S$ be a noetherian local scheme, $s$ the closed point
of $S$, and $f:X\to S$ a proper morphism of schemes. Suppose
that $\dim X(s)\leq 1$, and $H^1(X(s),\cO_{\!X(s)})=0$. Then
the natural map
$$
\Pic\,X\to\Pic\,X(s)
$$
is injective.
\end{claim}
\begin{pfclaim} Say that $S=\Spec\,A$ for a local ring
$A$, and denote by $\fm_A\subset A$ the maximal ideal. For
every $k\in\N$, set $S_n:=\Spec\,A/\fm_A^{k+1}$, and let
$i_n:X_n:=X\times_SS_n\to X$ be the closed immersion. Let
$\cL$ be any invertible $\cO_{\!X}$-module, and suppose
that $i_0^*\cL\simeq\cO_{\!X_0}$; we have to show that
$\cL\simeq\cO_{\!X}$. We notice that, for every $k\in\N$,
the natural map
$$
H^0(X_{k+1},i^*_{k+1}\cL)\to H^0(X_k,i^*_k\cL)
$$
is surjective : indeed, its cokernel is an $A$-submodule of
$H^1(X_0,\fm^k_A\cL/\fm^{k+1}_A\cL)$, and since $\dim X_0\leq 1$,
the natural map
$$
(\fm^k_A/\fm_A^{k+1})\otimes_{\kappa(s)}H^1(X_0,i_0^*\cL)\isom
H^1(X_0,(\fm^k_A/\fm_A^{k+1})\otimes_{\kappa(s)}\cL)\to
H^1(X_0,\fm^k_A\cL/\fm^{k+1}_A\cL)
$$
is surjective; on the other hand, our assumptions imply that
$H^1(X_0,i_0^*\cL)=0$, whence the contention. Let $A^\wedge$
be the $\fm_A$-adic completion of $A$; taking into account
\cite[Ch.III, Th.4.1.5]{EGAIII}, we deduce that the natural
map
$$
H^0(X,\cL)\otimes_AA^\wedge\to H:=H^0(X_0,i^*\cL)
$$
is a continuous surjection, for the $\fm_A$-adic topologies.
Since $H$ is a discrete space for this topology, and since
the image of $H^0(X,\cL)$ is dense in the $\fm_A$-adic
topology of $H^0(X,\cL)\otimes_AA^\wedge$, we conclude that
the restriction map $H^0(X,\cL)\to H$ is surjective as well.
Let $\bar s\in H$ be a global section of $i_0^*\cL$ whose image
in $i^*_0\cL_x$ is a generator of the latter $A_0$-module, for
every $x\in X_0$, and pick $s\in H^0(X,\cL)$ whose image
in $H$ equals $\bar s$. It remains only to check that, for
every $x\in X$, the $A$-module $\cL_x$ is generated by the
image $s_x$ of $s$. However, since $X$ is proper, every
$x\in X$ specializes to a point of $X_0$, hence we may
assume that $x\in X_0$, in which case one concludes easily,
by appealing to Nakayama's lemma (details left to the reader).
\end{pfclaim}

Combining claim \ref{cl_from-thesis} and \cite[Prop.11.1(i)]{Lip},
we see that the induced map
$$
i^*:\Pic\,Y'\to\Pic\,Y'_0
$$
is injective. Let now $\cL$ be any invertible
$\cO_{\!U}$-module; we have to show that $\cL$ extends to an
invertible $\cO_Y$-module (which is then isomorphic to $\cO_Y$).
However, notice that $f$ restricts to an isomorphism
$g:U'\isom U$, hence $g^*\cL$ is an invertible
$\cO_{\!U'}$-module, and by the foregoing there exists an
invertible $\cO_{\!Y'}$-module $\cL'$ such that
$\cL'_{|U'}\simeq g^*\cL$. In light of the isomorphism
\eqref{eq_this-is-h}, there exists an invertible
$\cO_{\P^1_Y}$-module $\cL''$ such that
$i^*\cL'\simeq i^*h^*\cL''$. Therefore
$$
\cL'\simeq h^*\cL''.
$$
Now, on the one hand, claim \ref{cl_from-thesis} implies that
$\cL''=\cO_{\P_Y^1}(n)$ for some $n\in\N$; on the other hand,
\eqref{eq_trivial-proj-sp} implies that $h$ restricts to a
morphism of schemes $h_\tr:U'\to\G_{m,Y}$, and
$\cL''_{|\G_{m,Y}}=\cO_{\P_Y^1}(n)_{|\G_{m,Y}}=\cO_{\G_{m,Y}}$,
so finally $g^*\cL=\cO_{\!U'}$, hence $\cL=\cO_{\!U}$, whence
the contention. 
\end{proof}

The following result complements proposition
\ref{prop_log-struct-is-fixed}.

\begin{theorem}\label{th_comp-fixed}
Let $(X,\underline M)$ be any regular log scheme, set
$U:=:(X,\underline M)_\tr$ and denote by $j:U\to X$ the open
immersion. Then we have :
$$
R^1j_*\cO^\times_{\!U}=0.
$$
\end{theorem}
\begin{proof} We begin with the following general :

\begin{claim}\label{cl_open-irredu}
Let $\pi:T_1\to T_2$ be a continuous open and surjective map of
topological spaces, such that $\pi^{-1}(t)$ is an irreducible
topological space (with the subspace topology) for every
$t\in T_2$. Then, we have :
\begin{enumerate}
\item
For every sheaf $S$ on $T_2$, the natural map $S\to\pi_*\pi^*S$
is an isomorphism.
\item
For every abelian sheaf $S$ on $T_2$, the natural map
$S[0]\to R\pi_*\pi^*S$ is an isomorphism in $\sD(\Z_{T_2}\Mod)$.
\end{enumerate}
\end{claim}
\begin{pfclaim}(i): Since $\pi$ is open, $\pi^*S$ is the sheaf
associated to the presheaf : $U\mapsto S(\pi U)$, for every
open subset $U$ of $T_1$. We show, more precisely, that this
presheaf is already a sheaf; since $\pi$ is surjective, the
claim shall follow immediately.
Now, let $U\subset T_1$ be an open subset, and $(U_i~|~i\in I)$
a family of open subsets of $X$ covering $U$; for every $i,j\in I$,
set $U_{ij}:=U_i\cap U_j$. It suffices to show that $S(\pi U)$
is the equalizer of the two maps :
$$
\xymatrix{
\prod_{i\in I}S(\pi U_i) \ar@<.5ex>[r] \ar@<-.5ex>[r]  &
\prod_{i,j\in I}S(\pi U_{ij}).
}$$
Since $S$ is a sheaf, the latter will hold, provided we know
that $\pi U_i\cap\pi U_j=\pi U_{ij}$ for every $i,j\in I$.
Hence, let $t\in\pi U_i\cap\pi U_j$; this means that
$\pi^{-1}(t)\cap U_i\neq\emptyset$ and
$\pi^{-1}(t)\cap U_j\neq\emptyset$. Since $\pi^{-1}(t)$ is
irreducible, we deduce that $\pi^{-1}(t)\cap U_{ij}\neq\emptyset$,
as required.

(ii): The proof of (i) also shows that, for every flabby abelian
sheaf $J$ on $T_2$, the abelian sheaf $\pi^*J$ is flabby on $T_1$.
Hence, if $S$ is any abelian sheaf on $T_2$, we obtain a flabby
resolution of $\pi^*S$ of the form $\pi^*J_\bullet$, by taking
a flabby resolution $S\to J_\bullet$ of $S$ on $T_2$. According
to remark \ref{rem_acyclic-crit}(iv), there is a natural isomorphism
$$
\pi_*\pi^*J\isom R\pi_*\pi^*S
$$
in $\sD(\Z_{T_2}\Mod)$. Then the assertion follows from (i).
\end{pfclaim}

After these preliminaries, let us return to the log scheme
$(X,\underline M)$, and its associated object
$\underline X:=((X,\underline M),F(X,\underline M),\pi_X)$.
The assertion to prove is local on $X$, hence we may assume
that $F(X,\underline M)=(\Spec\,P)^\sharp$, for some sharp,
fine and saturated monoid $P$. Next, by theorem
\ref{th_resolution-Zar} (and its proof) there exists an
integral proper simplicial subdivision
$\phi:F'\to F(X,\underline M)$, such that the log scheme
$(X',\underline M')$ underlying $\phi^*\underline X$ is
regular, $X'$ is regular, and the morphism $X'\to X$ restricts
to an isomorphism $U':=(X',\underline M')_\tr\to U$. In this
situation, we may find a further subdivision $\phi':F''\to F'$
such that both $\phi'$ and $\phi\circ\phi'$ are compositions of
saturated blow up of ideals generated by at most two elements
of $P$ (example \ref{ex_cut-by-hyperplane}(iii)). Say that
$\phi\circ\phi'=\phi_r\circ\phi_{r-1}\cdots\circ\phi_1$, where
each $\phi_i$ is a saturated blow up of the above type.
By proposition \ref{prop_basis-technique}(iii), we deduce
a sequence of morphisms of log schemes
$$
(X_1,\underline M{}_1)\xrightarrow{\ g_1\ }\cdots\to
(X_{r-1},\underline M{}_{r-1})\xrightarrow{\ g_{r-1}\ }
(X_r,\underline M{}_r)\xrightarrow{\ g_r\ }
(X_{r+1},\underline M{}_{r+1}):=(X,\underline M)
$$
each of which is the blow up of a corresponding ideal,
and by the same token, $\phi'$ induces a morphism
$g:(X_1,\underline M{}_1)\to(X',\underline M')$ of
$(X,\underline M)$-schemes. For every $i=1,\cdots,r+1$,
set $U_i:=(X_i,\underline M{}_i)_\tr$, and let
$j_i:U_i\to X_i$ be the open immersion; especially,
$U_{r+1}=U$. We shall show, by induction on $i$, that
\set\begin{equation}\label{eq_by-indu}
R^1j_{i*}\cO^\times_{U_i}=0
\qquad
\text{for $i=1,\dots,r+1$}.
\end{equation}
Notice first that the stated vanishing translates the
following assertion. For every $x\in X_i$ and every
invertible $\cO_{\!U_i}$-module $\cL$, there exists an
open neighborhood $U_x$ of $x$ in $X_i$ such that
$\cL_{|U_x\cap U_i}$ extends to an invertible
$\cO_{\!U_x}$-module. However, it follows immediately from
propositions \ref{prop_extend-rflx} and \ref{prop_Pic-is-Div},
that every invertible $\cO_{\!U'}$-module extends to an
invertible $\cO_{\!X'}$-module. Since $g$ restricts to
an isomorphism $U_1=g^{-1}U'\isom U'$, we easily deduce
that every invertible $\cO_{\!U_1}$-module extends to
an invertible $\cO_{\!X_1}$-module (namely : if $\cL$
is such a module, extend $g_{|U_1*}\cL$ to an invertible
$\cO_{\!X'}$-module, and pull the extension back to $X_1$,
via $g^*$). Summing up, we see that \eqref{eq_by-indu}
holds for $i=1$.

Next, suppose that \eqref{eq_by-indu} has already been
shown to hold for a given $i\leq r$; by lemma
\ref{lem_conclude-indu}(ii), it follows that
\eqref{eq_by-indu} holds for $i+1$, so we are done.
\end{proof}

\sset\subsubsection{}\label{subsec_from-P-to-X}
Let $X$ be a scheme, $\underline M$ a fine log structure on the
Zariski site of $X$, $x$ a point of $X$, and notice that every
fractional ideal of $\underline M{}_x$ is finitely generated (lemma
\ref{lem_rflx_and_quot}(iv)). Say that $X(x)=\Spec\,A$; in view of
lemma \ref{lem_rflx-rflx}(ii) there is a natural map of abelian
groups
$$
\Div(\underline M{}_x)\to\Div(A).
$$
Composing with the map $I\mapsto I^\sim$ as in
\eqref{eq_two-diff-rflx}, we get, by virtue of {\em loc.cit.},
a map
\set\begin{equation}\label{eq_from-P-to-X}
\Div(\underline M{}_x)\to\Div\,X(x).
\end{equation}

\begin{corollary} With the notation of \eqref{subsec_from-P-to-X},
suppose furthermore that $(X,\underline M)$ is regular at the
point $x$. Then \eqref{eq_from-P-to-X} induces an isomorphism
$$
\bar\Div(\underline M{}_x)\isom\bar\Div(X(x)).
$$
\end{corollary}
\begin{proof} By virtue of proposition \ref{prop_I-I}(ii) (and
remark \ref{rem_classy-rflx}(i)), the map under investigation
is already known to be injective. To show surjectivity, let
$K$ be the field of fractions of $A$, and $L\subset K$ any
reflexive fractional ideal of $A$; we may then regard
$\cL:=L^\sim$ as a coherent $\cO_{\!X}$-submodule of
$i_*\cO_{\!X_0}$, where $i:X_0\to X$ is the inclusion map
of the subscheme $X_0:=\Spec\,K$. By proposition
\ref{prop_Pic-is-Div}, the $\cO_{\!U}$-module $\cL_{|U}$ is
invertible, hence it extends to an invertible
$\cO_{\!X(x)}$-module $\cL'$, by virtue of theorem
\ref{th_comp-fixed}. Since $X(x)$ is local, $\cL'$ is a free
$\cO_{\!X(x)}$-module, and therefore $\cL_{|U}\simeq\cO_{\!U}$.
Thus, pick $a\in\cL(U)\subset K$ which generates $\cL_{|U}$;
after replacing $L$ by $a^{-1}L$, we may assume that
$\cL_{|U}=\cO_{\!U}$ as subsheaves of $i_*\cO_{\!X_0}$.
Let $\Sigma$ be the set of points of height one of $X(x)$
contained in $X(x)\setminus U$; for each $y\in\Sigma$,
the maximal ideal $\fm_y$ of $\cO_{\!X(x),y}$ is generated
by a single element $a_y$, and there exists $k_y\in\Z$ such
that $a_y^{k_y}$ generates the $\cO_{\!X(x),y}$-submodule
$\cL_y$ of $K$. To ease notation, let $P:=\underline M{}_x$,
and denote $\psi:X(x)\to\Spec\,P$ the natural continuous map;
also, let $\fm_{\psi(y)}$ be the maximal ideal of the localization
$P_{\psi(y)}$ for every $y\in X(x)$, and set
$$
I:=\bigcap_{y\in\Sigma}\fm_{\psi(y)}^{k_y}.
$$
In light of lemma \ref{lem_rflx_and_quot}(i,ii) and proposition
\ref{prop_classify-reflex-mon}, it is easily seen that
$I/P^\times$ is a reflexive fractional ideal of the fine and
saturated monoid $P^\sharp$, so $I$ is a reflexive fractional
ideal of $P$. Then $IA\subset K$ is a reflexive fractional ideal
of $A$. Set $\cL'':=(IA)^\sim\subset i_*\cO_{\!X_0}$; then
$\cL''_{|U}=\cL_{|U}$ and $\cL''_y=\cL_y$ for every $y\in\Sigma$.
It follows that, for every $y\in\Sigma$, there exists an open
neighborhood $U_y$ of $y$ in $X(x)$ such that
$\cL''_{|U_y}=\cL_{|U_y}$. Let $U':=U\cup\bigcup_{y\in\Sigma}U_y$,
and denote by $j':U'\to X(x)$ the open immersion. Notice that
$\delta'(y,\cO_{\!X})\geq 2$ for every $y\in X\setminus U'$
(corollary \ref{cor_normal-and-CM}). In light of proposition
\ref{prop_extend-rflx}(ii) (and remark \ref{rem_classy-rflx}(i)),
we deduce
$$
\cL=j'_*\cL_{|U'}=j'_*\cL''_{|U'}=\cL''
$$
whence the contention.
\end{proof}

\sset\subsubsection{}\label{subsec_not-yet-excel}
Let $\underline M$ be a log structure on the Zariski site
of a local scheme $X$, such that $(X,\underline M)$ is a regular
log scheme. Let $x\in X$ be the closed point, say that
$X=\Spec\,B$ for some local ring $B$, and let $\beta:P\to B$
a chart for $\underline M$ which is sharp at $x$.
As usual, if $M$ is any $B$-module, we denote $M^\sim$ the
quasi-coherent $\cO_{\!X}$-module arising from $M$.

\begin{theorem}
In the situation of \eqref{subsec_not-yet-excel}, suppose
as well that $\dim P=2$ and $\dim X=2$. Then every
indecomposable reflexive $\cO_{\!X}$-module is isomorphic
to $(IB)^\sim$, for some reflexive fractional ideal $I$ of $P$.
(Notation of \eqref{subsec_frac-ideal-mons}.)
\end{theorem}
\begin{proof} Set $Q:=\Div_+(P)$ and define $\phi:P\to Q$
as in example \ref{ex_Kummer-dim-two}. The chart $\beta$
defines a morphism $\psi:(X,\underline M)\to\Spec(\Z,P)$,
and we let $(X',\underline M')$ be the fibre product in the
cartesian diagram
$$
\xymatrix{ (X',\underline M') \ar[r] \ar[d]_f &
\Spec(\Z,Q) \ar[d]^{\Spec(\Z,\phi)} \\
(X,\underline M) \ar[r]^-\psi & \Spec(\Z,P).
}$$
Arguing as in the proof of claim \ref{cl_Kummer-preserve-regul},
it is easily seen that $X'$ is a local scheme and $(X',\underline M')$
is regular. Moreover, $X'$ is a regular scheme, since $Q$ is a
free monoid (corollary \ref{cor_more-precisely}), and $f$ is a
finite morphism of Kummer type (lemma \ref{lem_needed_once}).

\begin{claim}\label{cl_sogni-doro}
The $\cO_{\!X}$-module $f_*\cO_{\!X'}$ is isomorphic to a finite
direct sum $(I_1B\oplus\cdots\oplus I_kB)^\sim$, where
$I_1,\dots,I_k$ are reflexive fractional ideals of $P$, and
$\cO_{\!X}$ is a direct summand of $f_*\cO_{\!X'}$.
\end{claim}
\begin{pfclaim} In view of example \ref{ex_Kummer-dim-two}, we see
that $f_*\cO_{\!X'}=(Q\otimes_PB)^\sim$ is the direct sum of
the $\cO_{\!X}$-modules $(\gr_\gamma Q\otimes_PB)^\sim$, where
$\gamma$ ranges over the elements of $Q^\gp/P^\gp$, and
$\gr_\bullet Q$ denotes the $\phi$-grading of $Q$. But for each
such $\gamma$, the natural map
$\gr_\gamma Q\otimes_PB\to\gr_\gamma Q\cdot B$ is an isomorphism,
since $B$ is $P$-flat (proposition \ref{prop_second-crit}(ii)).
This shows the first assertion, and the second is clear as well,
since $\gr_0Q=P$.
\end{pfclaim}

Denote by $x'$ the closed point of $X'$; from corollary
\ref{cor_same-height}(i) we see that $U:=(X,\underline M)_1$ is
the complement of $\{x\}$ and $U':=(X',\underline M')_1$ is the
complement of $\{x'\}$ (notation of definition
\ref{def_trivial-locus}(i)).
Let also $j:U\to X$ and $j':U'\to X'$ be the open immersions.

\begin{claim}\label{cl_flip-page}
(i)\ \
The restriction $g:U'\to U$ of $f$ is a flat morphism of schemes.
\begin{enumerate}
\addenu
\item
The functor $j^*:\cO_{\!X}\bRflx\to\cO_{\!U}\bRflx$ is an
equivalence (see definition \ref{def_reflex}(iii)).
\end{enumerate}
\end{claim}
\begin{pfclaim}(i): It suffices to check that the restriction
$\Spec(\Z,Q)_1\to\Spec(\Z,P)_1$ of $\Spec(\Z,\phi)$ is flat.
However, we have the affine open covering
$$
\Spec(\Z,P)_1=\Spec\,\Z[P_{\fp_1}]\cup\Spec\,\Z[P_{\fp_2}]
$$
where $\fp_1,\fp_2\subset P$ are the two prime ideals of
height one (see example \ref{ex_satu-dim-two}(i)). Hence, we
are reduced to showing that the morphism of log schemes
underlying
$$
\Spec(\Z,\phi_{\fp_i}):\Spec(\Z,Q_{\fp_i})\to\Spec(\Z,P_{\fp_i})
$$
is flat for $i=1,2$. However, it is clear that $Q_{\fp_i}$
is an integral $P_{\fp_i}$-module, hence it suffices to check
that $Q_{\fp_i}$ is a flat $P_{\fp_i}$-module, for $i=1,2$
(proposition \ref{prop_Lazard}), or equivalently, that
$Q_{\fp_i}^\sharp$ is a flat $P_{\fp_i}^\sharp$-module
(corollary \ref{cor_yet-another-flat}(ii)). The latter
assertion follows immediately from the discussion of
\eqref{subsec_ramif-index}.

(ii): From proposition \ref{prop_extend-rflx} we see that
$j^*$ is full and essentially surjective. Moreover, it
follows from remark \ref{rem_co-repres} that every reflexive
$\cO_{\!X}$-module is $S_1$, so $j^*$ is also faithful
(details left to the reader).
\end{pfclaim}

In light of claim \ref{cl_flip-page}(ii), it suffices to show
that every indecomposable reflexive $\cO_{\!U}$-module $\cF$
is isomorphic to $(IB)^\sim_{|U}$, for some reflexive fractional
ideal $I$ of $P$. However, for such $\cF$, claim
\ref{cl_flip-page}(i) and lemma \ref{lem_rflx-on-lim}(i)
imply that $g^*\cF$ is a reflexive $\cO_{\!U'}$-module.
From proposition \ref{prop_extend-rflx} and corollary
\ref{cor_local-depth} we deduce that $\delta'(x',j'_*g^*\cF)\geq 2$,
so $j'_*g^*\cF$ is a free $\cO_{\!X'}$-module of finite rank
(\cite[Ch.0, Prop.17.3.4]{EGAIV}), and finally, $g^*\cF$ is a
free $\cO_{\!U'}$-module of finite rank. Taking into account
claim \ref{cl_sogni-doro}, it follows that
$g_*g^*\cF=\cF\otimes_{\cO_{\!U}}g_*\cO_{\!U'}$ is a direct
sum of $\cO_{\!U}$-modules of the type $(IB)^\sim_{|U}$, for
various $I\in\Div(P)$; moreover, $\cF$ is a direct summand of
$g_*g^*\cF$. Then, we may find a decomposition
$g_*g^*\cF=\cF_1\oplus\cdots\oplus\cF_t$ such that $\cF_i$ is
an indecomposable $\cO_{\!U}$-module for $i=1,\dots,t$, and
$\cF_1=\cF$ (details left to the reader : notice that -- since
reflexive $\cO_{\!U}$-modules are $S_1$ -- the length $t$ of
such a decomposition is bounded by the dimension of the
$\kappa(\eta)$-vector space $(g_*g^*\cF)_\eta$, where $\eta$
is the maximal point of $X$). On the other hand, notice that
$$
\begin{aligned}
\End_{\cO_{\!U}}((IB)^\sim_{|U})=\: & \End_B(IB) & &
\qquad\text{(claim \ref{cl_flip-page}(ii))} \\
=\: & (IB:IB) & &
\qquad\text{(by \eqref{subsec_J_P-plus})} \\
=\: & (I:I)B & &
\qquad\text{(lemma \ref{lem_rflx-rflx}(i))} \\
=\: & B & &
\qquad\text{(proposition \ref{prop_I-I}(i))}
\end{aligned}
$$
for every $I\in\Div(P)$. Now the contention follows from
theorem \ref{th_KRS}.
\end{proof}

\sset\subsubsection{}\label{subsec_down-from-log}
Let now $(X_\et,\underline M)$ be a quasi-coherent log
scheme on the \'etale site of $X$. Pick a covering family
$(U_\lambda~|~\lambda\in\Lambda)$ of $X$ in $X_\et$, and
for every $\lambda\in\Lambda$, a chart
$P_{\lambda,U_\lambda}\to\underline M_{|U_\lambda}$.
The latter induce isomorphisms of log schemes
$(U_\lambda,\underline M_{|U_\lambda})\isom
\tilde u{}^*(U_{\lambda,\Zar},P^{\log}_{U_{\lambda,\Zar}})$
(see \eqref{subsec_choose-a-top}); in other words, every
quasi-coherent log structure on $X_\et$ descends, locally on
$X_\et$, to a log structure on the Zariski site. However,
$(X_\et,\underline M)$ may fail to descend to a unique log
structure on the whole of $X_\Zar$. Our present aim is to
show that, at least under a few more assumptions on
$\underline M$, we may find a blow up
$(X'_\et,\underline M')\to(X_\et,\underline M)$ such that
$(X'_\et,\underline M')$ descends to a log structure on $X'_\Zar$.
To begin with, for every $\lambda\in\Lambda$ let
$$
T_\lambda:=(\Spec\,P_\lambda)^\sharp
\qquad\text{and}\qquad
S_\lambda:=\Spec\,\Z[P_\lambda].
$$
Also, let
$\underline S{}_\lambda:=
(\Spec(\Z,P_\lambda),T_\lambda,\psi_{P_\lambda})$ be the
object of $\cK$ attached to $P_\lambda$, as in
example \ref{ex_fan-andlogscheme}(i). From the isomorphism
$$
(U_{\lambda,\Zar},\tilde u_*\underline M_{|U_\lambda})\isom
U_\lambda\times_{S_\lambda}\Spec(\Z,P_\lambda)
$$
we deduce an object
$$
\underline U_\lambda:=
U_\lambda\times_{S_\lambda}\underline S{}_\lambda=
((U_{\lambda,\Zar},\tilde u_*\underline M_{|U_\lambda}),
T_\lambda,\psi_\lambda).
$$
Suppose now that $\underline M$ is a fs log structure; then
we may choose for each $P_\lambda$ a fine and saturated monoid
(lemma \ref{lem_simple-charts-top}(iii)).
Next, suppose that $X$ is quasi-compact; in this case we may
assume that $\Lambda$ is a finite set, hence
$\cS:=\{P_\lambda~|~\lambda\in\Lambda\}$ is a finite set of
monoids, and consequently we may choose a finite sequence of
integers $\underline c(\cS)$ fulfilling the conditions of
\eqref{subsec_models-of-fans} relative to the category
$\cS\text{-}\Fan$. Then, for every $\lambda\in\Lambda$ we have
a well defined integral roof $\rho_\lambda:T_\lambda(\Q_+)\to\Q_+$,
and we denote by $f_\lambda:T(\rho_\lambda)\to T_\lambda$ the
corresponding subdivision. There follows a cartesian morphism
$$
f^*_\lambda\underline U_\lambda\to\underline U_\lambda
$$
whose underlying morphism of log schemes is a saturated
blow up of the ideal
$\cI_{\!\rho_\lambda}\tilde u_*\underline M_{|U_\lambda}$
(proposition \ref{prop_to-subdiv-is-to-blow} and corollary
\ref{cor_case-of-blow-sat}). Next, for $\lambda,\mu\in\Lambda$,
set $U_{\lambda\mu}:=U_\lambda\times_X U_\mu$ and let
$$
\underline U_{\lambda\mu}:=
((U_{\lambda\mu},\tilde u_*\underline M_{|U_\lambda\mu}),
T_\lambda,\psi_{\lambda\mu})
$$
where $\psi_{\lambda\mu}$ is the composition of $\psi_\lambda$
and the projection $U_{\lambda\mu}\to U_\lambda$. Denote
by $(U^\sim_\lambda,\underline M^\sim_\lambda)$ (resp.
$(U^\sim_{\lambda\mu},\underline M^\sim_{\lambda\mu})$) the log
scheme underlying $f^*_\lambda\underline U_\lambda$ (resp.
$f_\lambda^*\underline U_{\lambda\mu}$). Also, for any
$\lambda,\mu,\gamma\in\Lambda$, let
$\pi_{\lambda\mu\gamma}:
U_{\lambda\mu}\times_X U_\gamma\to U_{\lambda\mu}$ be the natural
projection.

\begin{lemma}\label{lem_descent-for-tilde}
In the situation of \eqref{subsec_down-from-log}, we have :
\begin{enumerate}
\item
There exists a unique isomorphism of log schemes
$g_{\lambda\mu}:(U^\sim_{\lambda\mu},\underline M^\sim_{\lambda\mu})
\isom(U^\sim_{\mu\lambda},\underline M^\sim_{\mu\lambda})$
fitting into a commutative diagram :
$$
\xymatrix{ (U^\sim_{\lambda\mu},\underline M^\sim_{\lambda\mu})
           \ar[r]^-{g_{\lambda\mu}} \ar[d] &
           (U^\sim_{\mu\lambda},\underline M^\sim_{\mu\lambda})
           \ar[d] \\
           U_{\lambda\mu} \rdouble & U_{\mu\lambda}
}$$
whose vertical arrows are the saturated blow up morphisms.
\item
There exist natural isomorphisms of $\cO_{U_{\lambda\mu}}$-modules
$$
\omega_{\lambda\mu}:g^*_{\lambda\mu}\cO_{U_{\mu\lambda}}(1)\isom
\cO_{U_{\lambda\mu}}(1)
$$
such that \ \
$(\pi^*_{\mu\gamma\lambda}\omega_{\mu\lambda})\circ
(\pi^*_{\lambda\mu\gamma}\omega_{\lambda\mu})=
\pi^*_{\lambda\gamma\mu}\omega_{\lambda\gamma}$
\ \ for every $\lambda,\mu,\gamma\in\Lambda$.
\end{enumerate}
\end{lemma}
\begin{proof}(i): By the universal property of the saturated
blow up, it suffices to show that :
\set\begin{equation}\label{eq_coincide-here}
\cI_{\rho_\lambda}\tilde u_*\underline M_{|U_{\lambda\mu}}=
\cI_{\rho_\mu}\tilde u_*\underline M_{|U_{\mu\lambda}}
\qquad
\text{on\ \ $U_{\lambda\mu}=U_{\mu\lambda}$.}
\end{equation}
The assertion is local on $U_{\lambda\mu}$, hence let
$x\in U_{\lambda\mu}$ be any point; we get an isomorphism of
$\tilde u_*\underline M_x$-monoids :
$$
\cO_{T_\lambda,\psi_{\lambda\mu}(x)}\isom
\cO_{T_\mu,\psi_{\mu\lambda}(x)}
$$
whence an isomorphism of fans
$U(\psi_{\lambda\mu}(x))\isom U(\psi_{\mu\lambda}(x))$
(notation of \eqref{subsec_height-in-T}). Therefore, the subset
$U(x):=\psi_{\lambda\mu}^{-1}U(\psi_{\lambda\mu}(x))\cap
\psi^{-1}_{\mu\lambda}U(\psi_{\mu\lambda}(x))$
is an open neighborhood of $x$ in $U_{\lambda\mu}$, and both
$\psi_{\lambda\mu}$ and $\psi_{\mu\lambda}$ factor through
the same morphism of monoidal spaces :
$$
\psi(x):(U(x),(\tilde u_*\underline M^\sharp)_{|U(x)})
\to F(x):=(\Spec\,\tilde u_*\underline M_x)^\sharp
$$
and open immersions $F(x)\to T_\lambda$ and $F(x)\to T_\mu$.
It then follows from \eqref{subsec_models-of-fans} that
the preimage of $\cI_{\!\rho_\lambda}$ on $F(x)$ agrees
with the preimage of $\cI_{\!\rho_\mu}$, whence the contention.

(ii): By inspecting the definitions, it is easily seen that
the epimorphism \eqref{eq_especiallly} identifies naturally
$\cO_{U^\sim_{\lambda\mu}}(1)$ to the ideal
$\cI_{\!\rho_{\lambda\mu}}\cO_{U^\sim_{\lambda\mu}}$ of
$\cO_{U^\sim_{\lambda\mu}}$ generated by the image of
$\cI_{\rho_\lambda}\tilde u_*\underline M_{|U_{\lambda\mu}}$.
Hence the assertion follows again from \eqref{eq_coincide-here}.
\end{proof}

\sset\subsubsection{}\label{subsec_construct-tilde-X}
Lemma \ref{lem_descent-for-tilde} implies that
$$
((U^\sim_\lambda,\cO_{\!U^\sim_\lambda}(1)),
(g_{\lambda\mu},\omega_{\lambda\mu})~|~
\lambda,\mu\in\Lambda)
$$
is a descent datum -- relative to the faithfully flat and quasi-compact
morphism $\coprod_{\lambda\in\Lambda}U_\lambda\to X$ -- for the
fibred category of schemes endowed with an ample invertible sheaf.
According to \cite[Exp.VIII, Prop.7.8]{SGA1}, such a datum is effective,
hence it yields a projective morphism $\pi:X^\sim\to X$ together with
an ample invertible sheaf $\cO_{\!X^\sim}(1)$ on $X^\sim$, with
isomorphisms
$g_\lambda:X^\sim\times_X U_\lambda\isom U^\sim_\lambda$
of $U_\lambda$-schemes and $\pi_\lambda^*\cO_{\!X^\sim}(1)\isom
g^*_\lambda\cO_{\!U_\lambda^\sim}(1)$ of invertible modules.

Then the datum
$(\tilde u{}^*\underline M^\sim_\lambda,
\tilde u{}^*\log g_{\lambda\mu}~|~\lambda,\mu\in\Lambda)$
likewise determines a unique sheaf of monoids $\underline M^\sim$
on $X^\sim_\et$, and the structure maps of the log structures
$\underline M_\lambda$ glue to a well defined morphism of sheaves
of monoids $\underline M^\sim\to\cO_{\!X^\sim_\et}$, so that
$(X^\sim_\et,\underline M^\sim)$ is a log scheme, and the projection
$\pi$ extends to a morphism of log schemes
$(\pi,\log\pi):(X^\sim_\et,\underline M^\sim)\to(X_\et,\underline M)$.

\begin{proposition}\label{subsec_descends-to-Zar}
In the situation of \eqref{subsec_construct-tilde-X},
the counit of adjunction
$$
\tilde u{}^*\tilde u{}_*(X^\sim_\et,\underline M^\sim)
\to(X^\sim_\et,\underline M^\sim)
$$
is an isomorphism.
\end{proposition}
\begin{proof} (This is the counit of the adjoint pair
$(\tilde u{}^*,\tilde u_*)$ of \eqref{subsec_up_and_down-log},
relating the categories of log structures on $X^\sim_\Zar$ and
$X^\sim_\et$.) Recall that there exist natural epimorphisms
$(\Q_+^{\oplus d})_{T(\rho_\lambda)}\to \cO_{T(\rho_\lambda),\Q}$
(see \eqref{subsec_canonical-coordinates}), which induce epimorphisms
of $U^\sim_{\lambda,\Zar}$-monoids
$$
\theta_\lambda:(\Q_+^{\oplus d})_{U^\sim_{\lambda,\Zar}}\to
(\underline M^\sim_\lambda)^\sharp_\Q \qquad \text{for every
$\lambda\in\Lambda$}.
$$
The compatibility with open immersions expressed by
\eqref{eq_compat-coordinates} implies that the system
of maps $(\tilde u{}^*\theta_\lambda~|~\lambda\in\Lambda)$
glues to a well defined epimorphism of $X^\sim_\et$-monoids :
$$
\theta:(\Q_+^{\oplus d})_{X^\sim_\et}\to(\underline
M^\sim)^\sharp_\Q.
$$
In view of lemma \ref{lem_Hilbert90}(ii), it follows that the counit
of adjunction :
$$
\tilde u{}^*\tilde u{}_*(\underline M^\sim)^\sharp_\Q\to(\underline
M^\sim)^\sharp_\Q
$$
is an isomorphism. By applying again lemma \ref{lem_Hilbert90}(ii)
to the monomorphism $(\underline M^\sim)^\sharp\to(\underline
M^\sim)^\sharp_\Q$, we deduce that also the counit
$$
\tilde u{}^*\tilde u{}_*(\underline M^\sim)^\sharp\to(\underline
M^\sim)^\sharp
$$
is an isomorphism. Then the assertion follows from proposition
\ref{prop_reduce-to-etale}(iii).
\end{proof}

\begin{corollary} Let $(X_\et,\underline M)$ be a quasi-compact
regular log scheme. Then there exists a smooth morphism of log
schemes $(X'_\et,\underline M')\to(X_\et,\underline M)$ whose
underlying morphism of schemes is proper and birational, and such
that $X'$ is a regular scheme. More precisely, $f$ restricts
to an isomorphism
$(X'_\et,\underline M')_\tr\to(X_\et,\underline M)_\tr$ on the
trivial loci.
\end{corollary}
\begin{proof} Given such $(X_\et,\underline M)$, we construct first
the morphism $\pi:(X^\sim_\et,\underline M^\sim)\to(X_\et,\underline
M)$ as in \eqref{subsec_construct-tilde-X}. Since
$\pi\times_X\one_{U_\lambda}$ is proper for every
$\lambda\in\Lambda$, if follows that $\pi$ is proper (\cite[Ch.IV,
Prop.2.7.1]{EGAIV-2}). Likewise, notice that each morphism
$U^\sim_\lambda\to U_\lambda$ induces an isomorphism
$(U^\sim_\lambda,\underline M^\sim_\lambda)_\tr\isom
(U_\lambda,\tilde u_*\underline M_{|U_\lambda})_\tr$ (remark
\ref{rem_basic-technique}). It follows easily that $\pi$ restricts
an isomorphism on the trivial loci. Hence, we may replace
$(X_\et,\underline M)$ by $(X^\sim_\et,\underline M^\sim_\et)$.
Then, by corollary \ref{cor_undercover}(ii) and proposition
\ref{subsec_descends-to-Zar} we are further reduced to showing that
there exists a proper morphism $\pi':(X'_\Zar,\underline M')\to
\tilde u_*(X^\sim_\et,\underline M^\sim)$ with $X'$ regular,
restricting to an isomorphism on the trivial loci. However, in light
of lemma \ref{lem_reduce-to-et} (and again, proposition
\ref{subsec_descends-to-Zar}), the sought $\pi'$ is provided by the
more precise theorem \ref{th_resolution-Zar}.
\end{proof}

\subsection{Local properties of the fibres of a smooth morphism}
\label{sec_back-to-flatness}
Resume the situation of example \ref{ex_fan-andlogscheme}(ii),
and to ease notation set $\phi:=(\Spec\,\lambda)^\sharp$,
and $(f,\log f):=\Spec(R,\lambda)$.
Suppose now, that $\lambda:P\to Q$ is an integral, local
and injective morphism of fine monoids.
Then $f:S'\to S$ is flat and finitely presented (theorem
\ref{th_flat-crit-for-mnds}). Moreover :

\begin{lemma}\label{lem_dim-intgr}
In the situation of \eqref{sec_back-to-flatness}, for every
point $s\in S$, the fibre $f^{-1}(s)$ is either empty, or else
it is pure-dimensional, of dimension
$$
\dim f^{-1}(s)=\dim Q-\dim P=\rk_\Z\,\Coker\,\lambda^\gp.
$$
\end{lemma}
\begin{proof} To begin with, notice that
$\lambda^{-1}Q^\times=P^\times$, whence the second stated
identity, in view of corollary \ref{cor_consequent}(i).
To prove the first stated identity, we easily reduce to the case
where $R$ is a field. Notice that the image of $f$ is an open
subset $U\subset S$ (\cite[Ch.IV, Th.2.4.6]{EGAIV-2}), especially
$U$ (resp. $S'$) is pure-dimensional of dimension $\rk_\Z P^\gp$
(resp. $\rk_\Z Q^\gp$) by claim \ref{cl_NullSt}(ii).
Hence, fix any closed point $s\in S$, and set $X:=f^{-1}(s)$.
From \cite[Ch.IV, Cor.6.1.2]{EGAIV-2} we deduce that, for every
closed point $s'\in X$, the Krull dimension of $\cO_{\!X,s'}$
equals $r:=\rk_\Z\Coker\,\lambda^\gp$. More precisely, say that
$Z\subset X$ is any irreducible component; we may
find a closed point $s'\in Z$ which does not lie on any other
irreducible component of $X$, and then the foregoing implies
that the dimension of $Z$ equals $r$, as stated.

Next, let $s\in U$ be any point, and denote $K$ the residue
field of $\cO_{\!U,s}$, and $\pi:R[P]\to K$ the natural map.
Let $y\in\Spec\,K[P]$ be the $K$-rational closed point such
that $a(y)=\pi(a)$ for every $a\in P$; then the image of
$y$ in $S$ equals $s$, and if we let $f_K:=\Spec\,K[\lambda]$,
we have an isomorphism of $K$-schemes
$f_K^{-1}(y)\isom f^{-1}(s)$. The foregoing shows that
$f_K^{-1}(y)$ is pure-dimensional of dimension $r$, hence
the same holds for $f^{-1}(s)$.
\end{proof}

Now, let us fix $s\in S$, such that $f^{-1}(s)\neq\emptyset$,
and suppose that $\psi_P(s)=\fm_P$ is the closed point of $T_P$.
For every $\fq\in\phi^{-1}(\fm_P)$, the closure $\overline{\{\fq\}}$
of $\{\fq\}$ in $T_Q$ is the image of the natural map
$\Spec\,Q/\fq\to T_Q$. We deduce natural isomorphisms of schemes :
$$
S_0:=\psi_P^{-1}\{\fm_\fp\}\isom\Spec\,R\La P/\fm_P\Ra
\qquad
S'_\fq:=\psi_Q^{-1}\overline{\{\fq\}}\isom\Spec\,R\La Q/\fq\Ra
$$
under which, the restriction $f_\fq:S'_\fq\to S_0$ of
$f$ is identified with $\Spec\,R\La\lambda_\fq\Ra$, where
$\lambda_\fq:P/\fm_P\to Q/\fq$ is induced by $\lambda$.
The latter is an integral and injective morphism as well
(corollary \ref{cor_flat-face}). Explicitly, set
$F:=Q\!\setminus\!\fq$, and let $\lambda_F:P^\times\to F$
be the restriction of $\lambda$; then $\lambda_\fq=(\lambda_F)_\circ$,
and lemma \ref{lem_dim-intgr} yields the identity :
\set\begin{equation}\label{eq_dim-of-face}
\dim f^{-1}_\fq(s)=\dim Q/\fq=\dim Q-\hgt\,\fq.
\end{equation}
Also, notice that $T_Q$ is a finite set under the current
assumptions (lemma \ref{lem_face}(iii)), and clearly
$$
f^{-1}(s)=\!\!\!\!\!\!
\bigcup_{\fq\in\Max(\phi^{-1}\fm_P)}\!\!\!\!\!\!
f^{-1}_\fq(s)
$$
where, for a topological space $T$, we denote by $\Max(T)$
the set of maximal points of $T$.
Therefore, for every irreducible component $Z$ of $f^{-1}(s)$
there must exist $\fq\in\Max(\phi^{-1}\fm_P)$ such that
$Z\subset f^{-1}_\fq(s)$, and especially,
\set\begin{equation}\label{eq_Conversely}
\dim f^{-1}(s)=\dim f_\fq^{-1}(s).
\end{equation}
Conversely, we claim that \eqref{eq_Conversely} holds for
every $\fq\in\Max(\phi^{-1}\fm_P)$. Indeed, suppose that
$\dim f^{-1}_\fq(s)<\dim f^{-1}(s)$ for one such $\fq$, and
let $Z$ be an irreducible component of $f^{-1}_\fq(s)$; let
also $Z'$ be an irreducible component of $f^{-1}(s)$ containing
$Z$. By the foregoing, there exists
$\fq'\in\Max(\phi^{-1}\fm_P)$ with
$Z'\subset f_{\fq'}^{-1}(s)$. Set $\fq'':=\fq\cup\fq'$;
then clearly $\fq''\in\phi^{-1}(\fm_P)$, and
$$
\overline{\{\fq''\}}=\overline{\{\fq\}}\cap\overline{\{\fq'\}}.
$$
Especially, $Z\subset f^{-1}_{\fq''}(s)$; however, it follows
from \eqref{eq_dim-of-face} that
$\dim f^{-1}_{\fq''}(s)<\dim f^{-1}_\fq(s)$ (since
$\hgt\,\fq<\hgt\,\fq''$); but this is absurd,
since $\dim Z=\dim f^{-1}_\fq(s)$ (lemma \ref{lem_dim-intgr}).
The same counting argument also shows that every maximal point
of $f^{-1}(s)$ gets mapped necessarily to a maximal point of
$\phi^{-1}(\fm_P)$; in other words, we have shown that $\psi_Q$
restricts to a surjective map :
$$
\Max(f^{-1}s)\to\Max(\phi^{-1}\fm_P).
$$
More precisely, let $s'\in f^{-1}(s)$ be a point such that
$\psi_Q(s')=\fm_Q$, the closed point of $T_Q$. Then clearly
$s'\in f_\fq^{-1}(s)$ for every $\fq\in\Max(\phi^{-1}\fm_P)$,
and it follows that the foregoing surjection restricts further
to a surjective map :
\set\begin{equation}\label{eq_surjective-Max}
f^{-1}_{s'}(s)\to\Max(\phi^{-1}\fm_P).
\end{equation}
On the other hand, since $\log\psi_Q$ is an isomorphism, we have
$$
(Q^{\log}_{S'})^\sharp_{\bar s{}'}=\cO_{T_Q,\fq}=Q^\sharp_\fq
\qquad
\text{for every $s'\in\psi_Q^{-1}(\fq)$}
$$
(where $\bar s{}'$ denotes any $\tau$-point of $S'$ localized
at $s'$); explicitly, if $F=Q\!\setminus\!\fq$, then
$Q^\sharp_\fq=Q/F$; likewise, we get
$(P^{\log}_S)^\sharp_{\bar s}=P/\phi^{-1}F$.
Taking into account \eqref{eq_dim-of-face}, we deduce :
$$
\dim f^{-1}(f(s'))-\dim f_\fq^{-1}(f(s'))=
\rk_\Z\,\Coker\,(\log f)_{\bar s{}'}^\gp
\qquad
\text{for every $s'\in\psi_Q^{-1}(\fq)$.}
$$
Thus, for every $n\in\N$, let
$$
U_n:=\{s'\in S'~|~\rk_\Z\,\Coker\,(\log f)_{\bar s{}'}^\gp=n\}.
$$
The foregoing implies that for every $s\in S$, the subset
$U_0\cap f^{-1}(s)$ is open and dense in $f^{-1}(s)$, and
$U_n\cap f^{-1}(s)$ is either empty, or else it is a subset
of pure codimension $n$ in $f^{-1}(s)$.

\sset\subsubsection{}\label{subsec_geom-reduced-fibr}
In the situation of \eqref{sec_back-to-flatness}, suppose
moreover that $\log f$ is saturated; notice that this condition
holds if and only if $\log\phi$ is saturated, if and only if
the same holds for $\lambda$ (lemma \ref{lem_little}(iii)).
Then, corollary \ref{cor_persist-integr}(ii) says that
$\Coker(\log f^\sharp)^\gp_{\bar s{}'}$ is torsion-free
for every $s'\in S'$; especially,
$\Coker(\log f^\sharp)^\gp_{\bar s{}'}$ vanishes for every
$s'\in U_0$. Then corollary \ref{cor_persist-integr}(i) implies
that $\log f^\sharp_{\bar s{}'}$ is an isomorphism for every
$s'\in U_0$, in which case the same holds for $\log f_{\bar s{}'}$
(lemma \ref{lem_check-iso}); in other words, $U_0$ is the strict
locus of $f$ (see definition \ref{def_trivial-locus}(ii)). 

\begin{lemma}\label{lem_notnilpo}
Let $K$ be an algebraically closed field of characteristic $p$,
and $\chi:P\to(K,\cdot)$ a local morphism of monoids. Let also
$\lambda:P\to Q$ be as in  \eqref{subsec_geom-reduced-fibr}. We have :
\begin{enumerate}
\item
The $K$-algebra $Q\otimes_PK$ is Cohen-Macaulay.
\item
If moreover, either $p=0$, or else the order of the torsion
subgroup of\/ $\Gamma:=\Coker\,\lambda^\gp$ is not divisible
by $p$, then $Q\otimes_PK$ is reduced ({\em i.e.} its nilradical
is trivial).
\end{enumerate}
\end{lemma}
\begin{proof}(i): Since $\chi$ is a local morphism, we have
$\chi^{-1}(0)=\fm_P$, and $\chi$ is determined by its restriction
$P^\times\to K^\times$, which is a homomorphism of abelian
groups. Notice that $K^\times$ is divisible, hence it is injective
in the category of abelian groups; since the unit of adjunction
$P\to P^\sat$ is local (lemma \ref{lem_exc-satura}(iii)), it follows
that $\chi$ extends to a local morphism $\chi':P^\sat\to K$.
Notice that $Q\otimes_PK=Q^\sat\otimes_{P^\sat}K$ (lemma
\ref{lem_little}(iv)), hence we may replace $P$ by $P^\sat$
and $Q$ by $Q^\sat$, which allows to assume from start that
$Q$ is saturated. In this case, by theorem \ref{th_Hochster}(i),
both $K[P]$ and $K[Q]$ are Cohen-Macaulay; since $K[\lambda]$
is flat, theorem \ref{th_depth-flat-basechange} and
\cite[Ch.IV, Cor.6.1.2]{EGAIV-2} imply that $Q\otimes_PK$ is
Cohen-Macaulay as well.

(ii): Let $Q=\bigoplus_{\gamma\in\Gamma}Q_\gamma$ be the
$\lambda$-grading of $Q$ (see remark \ref{rem_why-not}(iii)).
Under the current assumptions, $\lambda$ is exact (lemma
\ref{lem_persist-integr}(ii)), consequently $Q_0=P$ (remark
\ref{rem_why-not}(v)). Moreover, $Q_\gamma$ is a finitely
generated $P$-module, for every $\gamma\in\Gamma$ (corollary
\ref{cor_no-fibres-here}), hence either $Q_\gamma=\emptyset$,
or else $Q_\gamma$ is a free $P$-module of rank one, say
generated by an element $u_\gamma$ (remark \ref{rem_why-not}(iv)).
Furthermore, $Q^k_\gamma=Q_{k\gamma}$ for every integer $k>0$
and every $\gamma\in\Gamma$ (proposition \ref{prop_crit-saturated}).
Thus, whenever $Q_\gamma\neq\emptyset$, the element
$u_\gamma\otimes 1$ is a basis of the $K$-vector space
$Q_\gamma\otimes_PK$, and $(u_\gamma\otimes 1)^k$ is a basis
of $Q_{k\gamma}\otimes_PK$, for every $k>0$; especially,
$u_\gamma\otimes 1$ is not nilpotent. In view of proposition
\ref{prop_notnilpo}(ii), the contention follows.
\end{proof}

\sset\subsubsection{}\label{subsec_was-in-th-segue}
Let $f:(X,\underline M)\to(Y,\underline N)$ be a smooth and
log flat morphism of fine log schemes. For every $n\in\N$, let
$U(f,n)\subset X$ be the subset of all $x\in X$ such that
$\rk_\Z\Coker\,(\log f^\sharp_{\bar x})^\gp=n$, for any
$\tau$-point $\bar x$ localized at $x$. By lemma
\ref{lem_simple-charts}(iii), $U(f,n)$ is a locally closed
subset (resp. an open subset) of $X$, for every $n>0$ (resp.
for $n=0$).

\begin{theorem}\label{th_satura-smooth}
In the situation of \eqref{subsec_was-in-th-segue}, we have :
\begin{enumerate}
\item
$f$ is a flat morphism of schemes.
\item
For all $y\in Y$, every connected component of $f^{-1}(y)$,
is pure dimensional, and $f^{-1}(y)\cap U_n$ is either
empty, or else it has pure codimension $n$ in $f^{-1}(y)$,
for every $n\in\N$.
\item
Moreover, if $f$ is saturated, we have :
\begin{enumerate}
\item
The strict locus $\mathrm{Str}(f)$ of $f$ is open in $X$,
and $\mathrm{Str}(f)\cap f^{-1}(y)$ is a dense subset
of $f^{-1}(y)$, for every $y\in Y$.
\item
$f^{-1}(y)$ is geometrically reduced and Cohen-Macaulay,
for every $y\in Y$.
\item
$\Coker(\log f^\sharp_{\bar x})^\gp$ is a free abelian group
of finite rank, for every $\tau$-point $\bar x$ of $X$.
\end{enumerate}
\end{enumerate}
\end{theorem}
\begin{proof} Let $\xi$ be any $\tau$-point of $X$;
according to corollary \ref{cor_first-trick}, there exist
a neighborhood $V$ of $\xi$, and a fine chart
$P_V\to\underline N{}_{|V}$, such that $P^\gp$ is a free
abelian group of finite rank, and the induced morphism
of monoids $P\to\cO_{Y,\xi}$ is local. 

By lemma \ref{lem_give-it-a-name},
\cite[Ch.IV, Th.2.4.6, Prop.2.5.1]{EGAIV-2} and
\cite[Ch.IV, Prop.17.5.7]{EGA4}, it suffices to prove the
theorem for the morphism
$f\times_YV:(X,\underline M)\times_YV\to(V,\underline N_{|V})$.
Hence we may assume from start that $\underline N$ admits a
fine chart $\beta:P_Y\to\underline N$, such that $P^\gp$ is
torsion-free and the induced map $P\to\cO_{Y,\xi}$ is local.

By theorem \ref{th_charact-smoothness}, remark
\ref{rem_strict-locus}(i), \cite[Ch.IV, Th.2.4.6]{EGAIV-2} and
\cite[Ch.IV, Prop.17.5.7]{EGA4} (and again lemma
\ref{lem_give-it-a-name}) we may further assume that $f$ admits
a fine chart $(\beta,\omega_Q:Q_X\to\underline M,\lambda)$, such
that $\lambda$ is injective, the torsion subgroup of
$\Coker\,\lambda^\gp$ is a finite group whose order is invertible
in $\cO_{\!X}$, and the induced morphism of schemes
$X\to Y\times_{\Spec\,\Z[Q]}\Spec\,\Z[P]$ is \'etale.
Moreover, by theorem \ref{th_good-charts}(iii), we may
assume -- after replacing $Q$ by a localization, and $X$
by a neighborhood of $\xi$ in $X_\tau$ -- that the morphism
$\lambda:P\to Q$ is integral (resp. saturated, if $f$ is saturated),
and the morphism $Q\to\cO_{\!X,\xi}$ induced by $\omega_{P,\xi}$
is local, so $\lambda$ is local as well.

In this case, the same sort of reduction as in the foregoing
shows that, in order to prove (i) and (ii), it suffices to
consider a morphism $f$ as in \eqref{sec_back-to-flatness},
for which these assertions have already been established.
Likewise, in order to show (iii), it suffices to consider
a morphism $f$ as in \eqref{subsec_geom-reduced-fibr}. For
such a morphism, assertions (iii.a) and (iii.c) are already
known, and (iii.b) is an immediate consequence of lemma
\ref{lem_notnilpo}.
\end{proof}

Theorem \ref{th_satura-smooth}(iii) admits the following
partial converse :

\begin{proposition} Let $f:(X,\underline M)\to(Y,\underline N)$
be a smooth and log flat morphism of fs log schemes, such
that $f^{-1}(y)$ is geometrically reduced, for every $y\in Y$.
Then $f$ is saturated.
\end{proposition}
\begin{proof} Fix a $\tau$-point $\xi$ of $X$; the assertion
can be checked on stalks, hence we may assume that $\underline N$
admits a fine and saturated chart $\beta:P_Y\to\underline N$
(lemma \ref{lem_simple-charts-top}(iii)), such that the induced
morphism $\alpha_P:P\to\cO_{Y,f(\xi)}$ is local (claim
\ref{cl_better-chart}).
Then, by corollary \ref{cor_charact-smoothness}, we may find an
\'etale morphism $g:U\to X$ and a $\tau$-point $\xi'$ of $U$
with $g(\xi')=\xi$, such that the induced morphism of log
schemes $f_U:(U,g^*\underline M)\to(Y,\underline N)$ admits a
fine and saturated chart
$(\beta,\omega_Q:Q_U\to g^*\underline M,\lambda)$, where $\lambda$
is injective, and the induced ring homomorphism
\set\begin{equation}\label{eq_rompi}
Q\otimes_P\cO_{Y,f(\xi)}\to\cO_{U,\xi'}
\end{equation}
is \'etale. By \cite[Ch.IV, Prop.17.5.7]{EGA4}, the fibres
of $f_U$ are still geometrically reduced, hence we are reduced
to the case where $U=X$ and $\xi=\xi'$. Furthermore, after
replacing $Q$ by a localization, and $X$ by a neighborhood
of $\xi$, we may assume that the map $\alpha_Q:Q\to\cO_{\!X,\xi}$
induced by $\omega_{Q,\xi}$ is local and $\lambda$ is integral
(theorem \ref{th_good-charts}(iii) and lemma \ref{lem_exc-satura}(i)).
Lastly, let $K$ be the residue field of $\cO_{Y,f(\xi)}$; our
assumption on $f^{-1}(y)$ means that the ring
$A:=\cO_{\!X,\xi}\otimes_{\cO_{\!Y,f(\xi)}}K$ is reduced.

We shall apply the criterion of proposition
\ref{prop_crit-saturated}. Thus, let
$Q=\bigoplus_{\gamma\in\Gamma}Q_\gamma$ be the $\lambda$-grading
of $Q$ and notice as well that $\lambda$ is a local morphism
(since the same holds for $\alpha_P$ and $\alpha_Q$),
therefore it is exact (lemma \ref{lem_persist-integr}(ii));
consequently $Q_0=P$ (remark \ref{rem_why-not}(v)). Moreover,
$Q_\gamma$ is a finitely generated $P$-module, for every
$\gamma\in\Gamma$ (corollary \ref{cor_no-fibres-here}), hence
either $Q_\gamma=\emptyset$, or else $Q_\gamma$ is a free
$P$-module of rank one. 

We have to prove that $Q^k_\gamma=Q_{k\gamma}$ for every integer
$k>0$ and every $\gamma\in\Gamma$. In case $Q_\gamma=\emptyset$,
this is the assertion that $Q_{k\gamma}=\emptyset$ as well.
However, since $Q$ is saturated, the same holds for
$\Gamma':=(\lambda P)^{-1}Q/(\lambda P)^\gp$ (lemma
\ref{lem_exc-satura}(i,ii)), so it suffices to remark that
$\Gamma'\subset\Gamma$ is precisely the submonoid consisting
of all those $\gamma\in\Gamma$ such that $Q_\gamma\neq\emptyset$.

Therefore, fix a generator $u_\gamma$ for every $P$-module
$Q_\gamma\neq\emptyset$, and by way of contradiction, suppose
that there exist $\gamma\in\Gamma'$ and $k>0$ such that
$u_\gamma^k$ does not generate the $P$-module $Q_{k\gamma}$;
this means that there exists $a\in\fm_P$ such that
$u_\gamma^k=a\cdot u_{k\gamma}$. Now, notice that the induced
morphism of monoids $\beta:P\to K$ is local, especially $\beta(a)=0$,
therefore $(u_\gamma\otimes 1)^k=0$ in the $K$-algebra $Q\otimes_PK$.
Denote by $I\subset Q\otimes_PK$ the annihilator of
$u_\gamma\otimes 1$, and notice that, since \eqref{eq_rompi} is flat,
$IA$ is the annihilator of the image $u'$ of $u_\gamma\otimes 1$ in $A$.

On the other hand, it is easily seen that
$I$ is the graded ideal generated by
$(u_\mu\otimes 1~|~\mu\in\Gamma'')$ where $\Gamma''\subset\Gamma'$
is the subset of all $\mu$ such that $u_\gamma\cdot u_\mu$ is
not a generator of the $P$-module $Q_{\gamma+\mu}$. Clearly
$u_\mu\notin Q^\times$ for every $\mu\in\Gamma''$, therefore
the image of $u_\mu\otimes 1$ in $A$ lies in the maximal ideal.
Therefore $IA\neq A$, {\em i.e.} $u'$ is a non-zero nilpotent
element, a contradiction. 
\end{proof}

\sset\subsubsection{}\label{subsec_acyclic-log}
Let $(X,\underline M)$ be any log scheme, and $\bar x$ any
geometric point, localized at a point $x\in X$.
Suppose that $y\in X$ is a generization of $x$, and $\bar y$
a geometric point localized at $y$; then, arguing as in
\eqref{subsec_strict-special}, we may extend uniquely any strict
specialization morphism $X(\bar y)\to X(\bar x)$ to a morphism
of log schemes
$(X(\bar y),\underline M(\bar y))\to(X(\bar x),\underline M(\bar x))$
fitting into a commutative diagram
$$
\xymatrix{ (X(\bar y),\underline M(\bar y)) \ar[r] \ar[d] &
(X(y),\underline M(y)) \ar[d] \\
(X(\bar x),\underline M(\bar x)) \ar[r] & (X(x),\underline M(x))
}$$
whose right vertical arrow is induced by the natural isomorphism
$$
(X(y),\underline M(y))\isom(X(x),\underline M(x))\times_{X(x)}X(y).
$$
A simple inspection shows that the induced morphism
$$
\Gamma(X(\bar x),\underline M(\bar x))^\sharp\to
\Gamma(X(\bar y),\underline M(\bar y))^\sharp
$$
is naturally identified with the morphism
$\underline M{}^\sharp_{\bar x}\to\underline M{}^\sharp_{\bar y}$
obtained from the specialization map
$\underline M{}_{\bar x}\to\underline M{}_{\bar y}$.

\sset\subsubsection{}\label{subsec_tardi}
Let $g:(X,\underline M)\to(Y,\underline N)$ be a morphism of
log schemes, $\bar x$ any geometric point of $X$, and set
$\bar y:=g(\bar x)$. The log structures
$\underline M\to\cO_{\!X}$ and $\underline N\to\cO_{\!Y}$,
and the map $\log g_{\bar x}$ induce a commutative diagram of
continuous maps
\set\begin{equation}\label{eq_sleep}
{\diagram
X(\bar x)\ar[r]^-{g_{\bar x}} \ar[d]_{\psi_{\bar x}} &
Y(\bar y) \ar[d]^{\psi_{\bar y}} \\
\Spec\,\underline M{}_{\bar x} \ar[r]^-{\phi_{\bar x}} &
\Spec\,\underline N{}_{\bar y}
\enddiagram}
\end{equation}
(notation of \ref{subsec_strict-loc-of-schs}), and notice
that $\psi_{\bar x}$ (resp. $\psi_{\bar y}$) maps the closed
point of $X(\bar x)$ (resp. of $Y(\bar y)$) to the closed
point $t_{\bar x}\in\Spec\,\underline M{}_{\bar x}$ (resp.
$t_{\bar y}\in\Spec\,\underline N{}_{\bar y}$).

\begin{proposition}\label{prop_max-pts-fibre}
In the situation of \eqref{subsec_tardi}, suppose that
$g$ is a smooth morphism of fine log schemes, and moreover :
\begin{enumerate}
\alphaenu
\item
either $g$ is a saturated morphism
\item
or $(X,\underline M)$ is a fs log scheme.
\end{enumerate}
Then the following holds :
\begin{enumerate}
\item
The map $\psi_{\bar x}$ restricts to a bijection :
$$
\Max(g^{-1}_{\bar x}(\bar y))\isom
\Max(\phi^{-1}_{\bar x}(t_{\bar y})).
$$
\item
For every irreducible component $Z$ of $g^{-1}_{\bar x}(\bar y)$,
set
$$
(Z,\underline M(Z)):=
(X(\bar x),\underline M(\bar x))\times_{X(\bar x)}Z.
$$
Then the $\kappa(\bar y)$-log scheme $(Z,\underline M(Z)_\red)$
is geometrically pointed regular. (Notation of example
{\em\ref{ex_push-trivial}(iv)} and remark
{\em\ref{rem_pointed-reg}(ii)}.)
\end{enumerate}
\end{proposition}
\begin{proof} By corollary \ref{cor_first-trick}, there
exist a neighborhood $U\to Y$ of $\bar y$, and a fine chart
$\beta:P_U\to\underline N{}_{|U}$ such that $\beta_{\bar y}$
is a local morphism, and $P^\gp$ is torsion-free. Now, let
$$
g_{\bar x}:=(g_{\bar x},\log g_{\bar x}):
(X(\bar x),\underline M(\bar x))\to
(Y(\bar y),\underline N(\bar y))
$$
be the morphism of log schemes induced by $g$ (notation of
\ref{subsec_acyclic-log}); by theorem \ref{th_charact-smoothness},
we may find a fine chart for $g_{\bar x}$ of the type
$(i^*_{\bar y}\beta,\omega,\lambda:P\to Q)$, where $\lambda$
is injective, and the order of the torsion subgroup of
$\Coker\,\lambda^\gp$ is invertible in $\cO_{\!X,x}$.
Moreover, set $R:=\cO_{Y(\bar y),\bar y}$; then the induced map
$X(\bar x)\to\Spec\,Q\otimes_PR$ is ind-\'etale, and -- after
replacing $Q$ by some localization -- we may assume that
$\omega_{\bar x}:Q\to\cO_{\!X(\bar x),\bar x}$ is local
(claim \ref{cl_better-chart}), hence the same holds for
$\lambda$. Furthermore, under assumption (a) (resp. (b)), we may
also suppose that $Q$ is saturated, by lemmata \ref{lem_exc-satura}(ii)
and \ref{lem_simple-charts-top}(ii) (resp. that $\lambda$ is saturated,
by theorem \ref{th_good-charts}(iii)).

Let us now define $f:S'\to S$ as in \eqref{sec_back-to-flatness};
notice that $\omega_{\bar x}$ induces a closed immersion
$Y(\bar y)\to S$, and we have a natural identification of
$Y(\bar y)$-schemes :
$$
\Spec\,Q\otimes_PR=Y(\bar y)\times_SS'.
$$
Denote by $\bar s$ the image of $\bar y$ in $S$, and $\bar s{}'$
the image of $\bar x$ in $Y(\bar y)\times_SS'\subset S'$; there
follows an isomorphism of $Y(\bar y)$-schemes :
$$
X(\bar x)\isom Y(\bar y)\times_{S(\bar s)}S'(\bar s{}')
$$
(\cite[Ch.IV, Prop.18.8.10]{EGA4}). Moreover, our chart induces
isomorphisms :
$$
\Spec\,\underline M{}_{\bar x}\isom T_Q
\qquad
\Spec\,\underline N{}_{\bar y}\isom T_P
$$
which identify $\phi_{\bar x}$ to the map $\phi:T_Q\to T_P$
of \eqref{sec_back-to-flatness}.
In view of these identifications, we see that \eqref{eq_sleep}
is the restriction to the closed subset $X(\bar x)$, of the
analogous diagram :
$$
\xymatrix{
S'(\bar s{}') \ar[r]^-{f_{\bar s{}'}} \ar[d]_{\psi_{\bar s{}'}} & 
S(\bar s) \ar[d]^{\psi_{\bar s}} \\
T_Q \ar[r]^\phi & T_P.
}$$
We may thus assume that $X=S'$, $Y=S$ and $g=f$. Moreover,
let $s$ (resp. $s'$) be the support of $\bar s$ (resp. of
$\bar s{}'$); the morphism $\pi:S'(\bar s{}')\to S'(s')$
is flat, hence it restricts to a surjection
$$
\Max(f^{-1}_{\bar s{}'}\bar s)\to\Max(f^{-1}_{s'}s).
$$
In order to prove (i), it suffices then to show that the map
$\Max(f^{-1}_{\bar s{}'}\bar s)\to\Max(\phi^{-1}\fm_P)$, defined
as the composition of the foregoing map and the surjection
\eqref{eq_surjective-Max}, is injective. This boils down to
the assertion that, for every $\fq\in\Max(\phi^{-1}\fm_P)$ the
point $s'$ lies on a unique irreducible component of the fibre
$\pi^{-1}(f_{\fq,s'}^{-1}s)$.
However, let $\bar\beta:P\to\kappa(s)$ be the composition of
the chart $P\to R$ and the projection $R\to\kappa(s)$; then
$f_\fq^{-1}(s)=\Spec\,Q/\fq\otimes_P\kappa(s)$. Since $\pi$
is ind-\'etale, the assertion will follow from
\cite[Ch.IV, Prop.17.5.7]{EGA4} and corollary \ref{cor_normal-and-CM},
together with :

\begin{claim}\label{cl_normal-strata}
The log scheme
$W_\fq:=\Spec\La\Z,Q/\fq\Ra\times_{\Spec\,\Z[P]}\Spec\,\kappa(s)$
is geometrically pointed regular.
\end{claim}
\begin{pfclaim} To ease notation, set $F:=Q\setminus\fq$;
by assumption $\lambda^{-1}F=P^\times$, so that
$$
W_\fq=(W'_\fq)_\circ
\qquad\text{where}\qquad
W'_\fq:=
\Spec(\Z,F)\times_{\Spec(\Z,P^\times)}\Spec\,\kappa(s)
$$
(notation of \eqref{subsec_constant-log}; notice that the log
structures of $\Spec(\Z,P^\times)$ and $\Spec\,\kappa(s)$
are trivial). Moreover, let $\lambda_F:P^\times\to F$ be the
restriction of $\lambda$; then $\lambda^\gp_F$ is injective,
and
$$
\Coker\,\lambda_F^\gp\subset\Coker\,\lambda^\gp
$$
(corollary \ref{coro_satur-face}(i)), hence the order of the
torsion subgroup of $\Coker\lambda_F^\gp$ is invertible in
$\cO_{\!S,s}$. Moreover, if $\lambda$ is saturated, then the
same holds for $\lambda_F$ (corollary \ref{coro_satur-face}(ii)),
and it is easily seen that if $Q$ is saturated, the same holds
for $F$.
Consequently, the morphism $\Spec(\Z,\lambda_F)$ is smooth
(proposition \ref{prop_toric-smooth}).
The same then holds for the morphism $W'_\fq\to\Spec\,\kappa(s)$
obtained after base change of $\Spec(\Z,\lambda_F)$ along the
morphism $h:\Spec(\Z,P^\times)\to\Spec\,\kappa(s)$ induced by
$\bar\beta$ (proposition \ref{prop_sorite-smooth}(ii)).

Now we notice that under either of the assumptions (a) or (b),
$W'_\fq$ is a fs log scheme. Indeed, under assumption (b),
this follows by remarking that $\Spec\,\kappa(s)$ is trivially
a fs log scheme, and $\Spec(\Z,\lambda_F)$ is saturated. Under
assumption (a), $\Spec(\Z,F)$ is a fs log scheme, and it
suffices to observe that $h$ is a strict morphism. Lastly,
since the morphism $W'_\fq\to\Spec\,\kappa(s)$ is obviously
saturated, we apply corollary \ref{cor_smooth-preserve-reg}
to conclude.
\end{pfclaim}

(ii): In light of the foregoing, we see that, for any irreducible
component $Z$ of $g^{-1}_{\bar x}(\bar y)$, there exists a unique
$\fq(Z)\in\Max(\phi^{-1}\fm_P)$ such that $Z$ is isomorphic to
the strict henselization of $f^{-1}_{\fq(Z)}(s)$, at the geometric
point $\bar s{}'$.
Notice now that the log structure of $W_{\fq(Z)}$ is reduced,
by virtue of claim \ref{cl_normal-strata} and proposition
\ref{prop_log-struct-is-fixed}; then, a simple inspection
of the definitions shows that $(Z,\underline M(Z)_\red)$
is isomorphic to the strict henselization $W_{\fq(Z)}(\bar s{}')$.
Invoking again claim \ref{cl_normal-strata}, we deduce the contention.
\end{proof}

\sset\subsubsection{}\label{subsec_almost-tardi}
In the situation of \eqref{subsec_tardi}, suppose that
$Y$ is a normal scheme, and $(Y,\underline N)_\tr$ is a dense
subset of $Y$. Let $\bar\eta$ be a geometric point of $Y(\bar y)$,
localized at the generic point $\eta$, and
$$
(U_\fq~|~\fq\in\Spec\,\underline M{}_{\bar x})
$$
the logarithmic stratification of $(X(\bar x),\underline M(\bar x))$
(see \eqref{subsec_log-stratif}). Notice that
$\psi_{\bar x}(g^{-1}(\eta))$ lies in the preimage
$\Sigma\subset\Spec\,\underline M{}_{\bar x}$ of the maximal
point $\emptyset$ of $\Spec\,\underline N{}_{\bar y}$;
therefore, $g^{-1}_{\bar x}(\eta)$ is the union of the subsets
$U_\fq\times_Y|\eta|$, for all $\fq\in\Sigma$. 

\begin{proposition}\label{prop_generic-boundary}
In the situation of \eqref{subsec_almost-tardi}, suppose that $g$
is a smooth morphism of fine log schemes. Then the following holds :
\begin{enumerate}
\item
$X$ is a normal scheme.
\item
The scheme $g^{-1}_{\bar x}(\bar\eta)$ is normal and irreducible.
\item
For every $\fq\in\Sigma$, the $\kappa(\eta)$-scheme $U_\fq\times_Y|\eta|$
is non-empty, geometrically normal and geometrically irreducible,
of pure dimension $\dim g^{-1}_{\bar x}(\eta)-\hgt\,\fq$.
\item
Especially, set
$W:=(X(\bar x)\!\setminus\!U_\emptyset)\times_Y|\bar\eta|$;
then $\psi_{\bar x}$ induces a bijection :
$$
\Max(W)
\isom\{\fq\in\Sigma~|~\hgt\fq=1\}.
$$
\item
For every $w\in\Max(W)$, the stalk
$\cO_{\!g^{-1}_{\bar x}(\bar\eta),w}$ is a discrete valuation ring.
\end{enumerate}
\end{proposition}
\begin{proof} Set $R:=\cO_{Y(\bar y),\bar y}$; arguing as in the
proof of proposition \ref{prop_max-pts-fibre}, we may find :
\begin{itemize}
\item
a local, flat and saturated morphism $\lambda:P\to Q$ of fine
monoids, such that the order of the torsion subgroup of
$\Coker\,\lambda^\gp$ is invertible in $R$;
\item
local morphisms of monoids $P\to R$, $Q\to\cO_{\!X(\bar x),\bar x}$
which are charts for the log structures deduced from $\underline N$
and respectively $\underline M$, and such that the induced morphism
of $Y(\bar y)$-schemes $X(\bar x)\to\Spec\,Q\otimes_PR$ is ind-\'etale.
\end{itemize}
By \cite[Ch.IV, Prop.17.5.7]{EGA4}, we may then assume that
$(X,\underline M)$ (resp. $(Y,\underline N)$) is the scheme
$\Spec\,Q\otimes_PR$ (resp. $\Spec\,R$), endowed with the log
structure deduced from the natural map $Q\to Q\otimes_PR$
(resp. the chart $P\to R$), and $g$ is the natural projection.
Suppose first that $R$ is excellent, and let $R'$ be the normalization
of $R$ in a finite extension $K'$ of $\Frac(R)$; then $R'$ is also
strictly local and noetherian, and if $y'$ denotes the closed point
of $Y':=\Spec\,R'$, then the residue field extension
$\kappa(y)\subset\kappa(y')$ is purely
inseparable. Set
$$
(X',\underline M'):=(X,\underline M)\times_YY'
\qquad
(Y',\underline N'):=(Y,\underline N)\times_YY'
$$
and let $g':(X',\underline M')\to(Y',\underline N')$ be the
induced morphism of log schemes; it follows especially that
the restriction $g^{\prime-1}(y')\to g^{-1}(y)$ is a homeomorphism
on the underlying topological spaces. Hence, there is a
geometric point $\bar x{}'$ of $X'$, unique up to isomorphism,
whose image in $X$ agrees with $\bar x$, and we easily deduce
an isomorphism of $Y$-schemes (\cite[Ch.IV, Prop.18.8.10]{EGA4})
\set\begin{equation}\label{eq_same-strict-hens}
X'(\bar x{}')\isom X(\bar x)\times_YY'.
\end{equation}
Let $\eta'$ be the generic point of $Y'$; by assumption,
$\underline N'$ is trivial in a Zariski neighborhood of $\eta'$,
hence $(X',\underline N')\times_{Y'}|\eta'|$ is a fs log schemes
(since $g$ is saturated), and then the same log scheme is also
regular (theorem \ref{th_smooth-preserve-reg}), therefore
$g^{\prime-1}(\eta')$ is a normal scheme (corollary
\ref{cor_normal-and-CM}). On the other hand, $R'$ is a Krull
domain (\cite[Th.33.10]{Na}), and $g'$ is flat with reduced
fibres (theorem \ref{th_satura-smooth}(iii.b)), so $X'$ is a
noetherian normal scheme (lemma \ref{lem_Krullu}), consequently
the same holds for $X'(\bar x{}')$
(\cite[Ch.IV, Prop.18.8.12(i)]{EGA4}); especially, the latter
is irreducible, so the same holds for
$X'(\bar x{}')\times_{Y'}|\eta'|$. In view of
\eqref{eq_same-strict-hens}, it follows that $X(\bar x)\times_Y|\eta'|$
is also normal and irreducible. Since $K'$ is arbitrary, this
completes the proof of (i) and (ii), in this case.

Next, if $R$ is any normal ring, we may write $R$ as
the union of a filtered family $(R_i~|~i\in I)$ of excellent
normal local subrings (\cite[Ch.IV, (7.8.3)(ii,vi)]{EGAIV-2}).
For each $i\in I$, denote by $\bar y_i$ the image of $\bar y$
in $\Spec\,R_i$; then the strict henselization $R^\sh_i$ of
$R_i$ at $\bar y_i$ is also a subring of $R$, so we may replace
$R_i$ by $R^\sh_i$, which allows to assume that each $R_i$ is
strictly local, normal and noetherian (
\cite[Ch.IV, Prop.18.8.8(iv), Prop.18.8.12(i)]{EGA4}). 
Up to replacing $I$ by a cofinal subset, we may assume that
the image of $P$ lies in $R_i$, for every $i\in I$.
For each $i\in I$, set $X_i:=\Spec\,Q\otimes_PR_i$,
$Y_i:=\Spec\,R_i$, and endow $X_i$ (resp. $Y_i$) with the
log structure $\underline M{}_i$ (resp. $\underline N{}_i$)
deduced from the natural map $Q\to Q\otimes_PR_i$ (resp.
$P\to R_i$). There follows a system of morphisms of log schemes
$g_i:(X_i,\underline M{}_i)\to(Y_i,\underline N{}_i)$ for
every $i\in I$, whose limit is the morphism $g$. Moreover,
since $(Y,\underline N)_\tr$ is dense in $Y$, the image of
$P$ lies in $R\!\setminus\!\{0\}$, hence it lies in
$R_i\setminus\!\{0\}$ for every $i\in I$, and the latter
means that $(Y_i,\underline N{}_i)_\tr$ is dense in $Y_i$,
for every $i\in I$. Let $\bar\eta_i$ (resp. $\bar x_i$)
be the image of $\bar\eta$ (resp. of $\bar x$) in $Y_i$
(resp. in $X_i$); moreover, for each $i\in I$, let $x_i\in X_i$
be the support of $\bar x_i$. By the previous case, we know
that $g_{i,\bar x_i}^{-1}(\bar\eta_i)$ is normal and irreducible.
However, $X$ (resp. $g^{-1}_{\bar x}(\bar\eta)$) is the limit of
the system of schemes $(X_i~|~i\in I)$ (resp.
$(g_{i,\bar x_i}^{-1}(\bar\eta_i)~|~i\in I)$), so (i) and
(ii) follow. (Notice that the colimit of a filtered system
of integral (resp. normal) domains, is an integral (resp.
normal) domain : exercise for the reader.)

(iii): For every $\fq\in\Spec\,Q=\Spec\,\underline M{}_{\bar x}$,
set $X_\fq:=\Spec\,Q/\fq\otimes_PR$; since the chart
$Q\to\cO_{X(\bar x),\bar x}$ is local, it is clear that
$x\in X_\fq$ for every such $\fq$, and :
$$
X_\fq(\bar x)=\bigcup_{\fp\in\Spec\,Q/\fq}U_\fp.
$$
If $\phi(\fq)=\emptyset$, the induced
map $P\to Q/\fq$ is still flat (corollary \ref{cor_flat-face}),
hence the projection $X_\fq(\bar x)\to Y$ is a flat morphism of
schemes, especially
$g^{-1}_{\bar x}(\eta)\cap X_\fq(\bar x)\neq\emptyset$.
Furthermore, the subset $X_\fq(\bar x)$ is pure-dimensional, of
codimension $\hgt\,\fq$ in $g^{-1}_{\bar x}(\eta)$, by
\eqref{eq_dim-of-face} and \cite[Ch.IV, Cor.6.1.4]{EGAIV-2}.
It follows that $U_\fq$ is a dense open subset of $X_\fq(\bar x)$.
To conclude, it remains only to show that each $X_\fq(\bar x)$ is
geometrically normal and geometrically irreducible; however, let
$j_\fq:X_\fq\to X$ be the closed immersion; the induced morphism
of log schemes $g_\fq:(X_\fq,j_\fq^*\underline M)\to(Y,\underline N)$
is also smooth, hence the assertion follows from (ii).

(iv) is a straightforward consequence of (iii).

(v): Notice that $A:=\cO_{\!g_{\bar x}^{-1}(\bar\eta),w}$ is
ind-\'etale over the noetherian ring $Q\otimes_P\kappa(\bar\eta)$,
hence its strict henselization is noetherian, and then $A$
itself is noetherian (\cite[Ch.IV, Prop.18.8.8(iv)]{EGA4}).
Since $X$ is normal, and $w$ is a point of height one
in $g^{-1}_{\bar x}(\bar\eta)$, we conclude that $A$ is a
discrete valuation ring.
\end{proof}

\section{\'Etale coverings of schemes and log schemes}
\label{chap_etale-cov}
\subsection{Acyclic morphisms of schemes}\label{sec_etale-cov}
For any scheme $X$, we shall denote by :
$$
\bCov(X)
$$
the category whose objects are the finite \'etale morphisms
$E\to X$; the morphisms $(E\to X)\to(E'\to X)$ are the
$X$-morphisms of schemes $E\to E'$. By faithfully flat descent,
$\bCov(X)$ is naturally equivalent to the subcategory of
$X^\sim_\et$ consisting of all locally constant constructible
sheaves. If $f:X\to Y$ is any morphism of schemes, and
$\varphi:E\to Y$ is an object $\bCov(Y)$, then
$f^*\varphi:=\varphi\times_YX:E\times_YX\to X$ is an object
of $\bCov(X)$; more precisely, we have a fibration :
\set\begin{equation}\label{eq_fibred-cat-cov}
\bCov\to\Sch
\end{equation}
over the category of schemes, whose fibre, over any scheme
$X$, is the category $\bCov(X)$.

\begin{lemma}\label{lem_fibred-cat-cov}
Let $f$ be a morphism of schemes, and suppose that either
one of the following conditions holds :
\begin{enumerate}
\alphaenu
\item
$f$ is integral and surjective.
\item
$f$ is faithfully flat.
\item
$f$ is proper and surjective.
\end{enumerate}
Then $f$ is of universal $2$-descent for the fibred
category \eqref{eq_fibred-cat-cov}.
\end{lemma}
\begin{proof} This is \cite[Exp.VIII, Th.9.4]{SGA4-2}.
\end{proof}

\begin{lemma}\label{lem_Lefschetz-covs}
In the situation of definition {\em\ref{def_Lef}}, let
$U\subset X$ be any open subset with $Y\subset U$. If\/
$\Lef(X,Y)$ holds, the closed immersion $j:Y\to U$ induces
a fully faithful functor
$$
j^*:\bCov(U)\to\bCov(Y).
$$
\end{lemma}
\begin{proof} Let $\bCov(\fX)$ be the full subcategory of
$\cO_\fX\Alg_\lfft$ consisting of all finite \'etale
$\cO_\fX$-algebras (notation of lemma \ref{lem_was-moved}(ii));
the category $\bCov(U)$ is a full subcategory of
$\cO_U\Alg_\lfft$, so lemma \ref{lem_was-moved}(ii) already
implies that the functor $\bCov(U)\to\bCov(\fX)$ is fully
faithful, hence we are reduced to showing that the functor :
$$
\bCov(\fX)\to\bCov(Y)
\qquad
\cA\mapsto\Spec\,\cA\otimes_{\cO_\fX}\cO_Y
$$
is fully faithful. To this aim, let $\cI\subset\cO_{\!X}$ be the
ideal defining the closed immersion $Y\subset X$; consider the
direct system of schemes
$$
(Y_n:=\Spec\,\cO_{\!X}/\!\cI^{n+1}~|~n\in\N)
$$
and let $\bCov(Y_\bullet)$ be the category consisting
of all direct systems $(E_n~|~n\in\N)$ of schemes, such
that $E_n$ is finite {\'e}tale over $Y_n$, and such that
the transition maps $E_n\to E_{n+1}$ induce isomorphisms
$E_n\isom E_{n+1}\times_{Y_{n+1}}Y_n$ for every $n\in\N$.
The morphisms in $\bCov(Y_\bullet)$ are the morphisms
of direct systems of schemes. We have a natural fully
faithful functor :
$$
\bCov(\fX)\to\bCov(Y_\bullet)
\qquad\cA\mapsto(\Spec\,\cA/\!\cI^{n+1}\cA~|~n\in\N)
$$
(\cite[Ch.I, Cor.10.6.10(ii)]{EGAI}). Finally, the functor
$\bCov(Y_\bullet)\to\bCov(Y)$ given by the rule :
$(E_n~|~n\in\N)\to E_0$ is an equivalence, by
\cite[Ch.IV, Th.18.1.2]{EGA4}. The claim follows.
\end{proof}

\sset\subsubsection{}\label{subsec_filtered-co-cov}
Consider now a cofiltered family $\cS:=(S_\lambda~|~\lambda\in\Lambda)$
of affine schemes. Denote by $S$ the limit of $\cS$, and
suppose moreover that $\Lambda$ admits a final element
$0\in\Lambda$. Let $f_0:X_0\to S_0$ be a finitely presented
morphism of schemes, and set :
$$
X_\lambda:=X_0\times_{S_0}S_\lambda
\qquad
f_\lambda:=f_0\times_{S_0}S_\lambda:
X_\lambda\to S_\lambda
\qquad
\text{for every $\lambda\in\Lambda$}.
$$
Set as well $X:=X_0\times_{S_0}S$ and $f:=f_0\times_{S_0}S:X\to S$.
For every $\lambda\in\Lambda$, let $p_\lambda:S\to S_\lambda$
(resp. $p'_\lambda:X\to X_\lambda$) be the natural morphism.
The functors $p^{\prime*}_\lambda:\bCov(X_\lambda)\to\bCov(X)$
define a pseudo-cocone in the $2$-category $\bCat$, whence
a functor :
\set\begin{equation}\label{eq_pseudo-at-last}
\Pscolim{\lambda\in\Lambda}\bCov(X_\lambda)\to\bCov(X).
\end{equation}

\begin{lemma}\label{lem_go-to-pseudo-lim}
In the situation of \eqref{subsec_filtered-co-cov},
the functor \eqref{eq_pseudo-at-last} is an equivalence.
\end{lemma}
\begin{proof} It is a rephrasing of
\cite[Ch.IV, Th.8.8.2, Th.8.10.5]{EGAIV-3}
and \cite[Ch.IV, Prop.17.7.8(ii)]{EGA4}.
\end{proof}

\begin{lemma}\label{lem_replace}
Let $X$ be a scheme, $j:U\to X$ an open immersion with
dense image, and $f:X'\to X$ an integral surjective and
radicial morphism. The following holds :
\begin{enumerate}
\item
The morphism $f$ induces an equivalence of sites
$$
f^*:X'_\et\to X_\et
$$
and the functor $f^*:\bCov(X)\to\bCov(X')$ is an equivalence.
\item
The functor $j^*:\bCov(X)\to\bCov(U)$ is faithful.
\item
Suppose that $X$ is reduced and normal, and moreover :
\begin{enumerate}
\item
either $X$ has finitely many maximal points
\item
or else, $j$ is a quasi-compact open immersion.
\end{enumerate}
Then $j^*$ is fully faithful, and its essential image
consists of all the objects $\phi$ of $\bCov(U)$ such
that $\phi\times_XX(\bar x)$ lies in the essential
image of the pull-back functor :
$$
\bCov(X(\bar x))\to\bCov(X(\bar x)\times_XU)
$$
for every geometric point $\bar x$ of $X$. (Notation of
definition {\em\ref{def_strict-loc}(ii)}.)
\item
Furthermore, if $X$ is locally noetherian and regular, and
$X\!\setminus\! U$ has codimension $\geq 2$ in $X$, then $j^*$
is an equivalence.
\end{enumerate}
\end{lemma}
\begin{proof}(i) follows from \cite[Exp.VIII, Th.1.1]{SGA4-2}.

(ii): Indeed, let $\phi:E\to X$ and $\phi':E'\to X$ be any
two finite \'etale morphisms, and $f,g:E\to E'$ two morphisms
of $X$-schemes such that $f\times_XU=g\times_XU$; let
$\Delta_{E'}\to E'\times_XE'$ be the open and closed diagonal
immersion, $(f,g):E\to E'\times_XE'$ the morphism deduced
from $f$ and $g$, and set $D:=(f,g)^{-1}\Delta_{E'}$. Then $D$
is the largest open subset of $E'$ such that $f_{|D}=g_{|D}$;
on the other hand, by assumption $\phi^{-1}U\subset D$, and
since $\phi$ is an open map, $\phi^{-1}U$ is dense in $E$.
Lastly, $D$ is also a closed subset of $E$, so $D=E$, whence
the claim.

(iii): Choose a covering $X=\bigcup_{i\in I}V_i$ consisting of
affine open subsets, let
$$
X':=\coprod_{i\in I}V_i
\qquad
X'':=X'\times_XX'
$$
and denote by $g:X'\to X$ the induced morphism; set also
$j'_i:=j\times_XV_i$ for every $i\in I$. By lemma
\ref{lem_fibred-cat-cov}, $g$ is of universal $2$-descent
for the fibred category $\bCov$. On the other hand, the
induced open immersion $j'':U\times_XX''\to X''$ has dense
image, hence the corresponding functor $j''{}^*$ is faithful,
by (ii). By corollary \ref{cor_first-theor}(ii), the full
faithfulness of $j^*$ follows
from the full faithfulness of the  pull-back functor $j'{}^*$
corresponding to the open immersion  $j':=j\times_XX'$. The
latter holds if and only if the same holds for all the pull-back
functors $j^{\prime*}_i$. Hence, we may replace $X$ by $V_i$,
and assume from start that $X$ is affine, say $X=\Spec\,A$.
Let $E\to X$ and $E'\to X$ be two objects of $\bCov(X)$, and
$h:E\times_XU\to E'\times_XU$ a morphism, and write
$E=\Spec\,B$, $E'=\Spec\,B'$ for finite \'etale $A$-algebras
$B$ and $B'$; we have to check that $h$ extends
to a morphism $E\to E'$.

$\bullet$\ \
To this aim, consider first the case where $X$ has finitely
many maximal points $\eta_1,\dots,\eta_s$. Under the
current assumptions, $A$ is the product of $s$ domains,
and its total ring of fractions $\Frac\,A$ is the product
of fields $\kappa(\eta_1)\times\cdots\times\kappa(\eta_s)$.
The restrictions
$h_{\eta_i}:=h\times_UX(\eta_i):E(\eta_i)\to E'(\eta_i)$
induce a map of $\Frac\,A$-algebras
$$
h^\natural_\eta:=\prod^s_{i=1}h^\natural_{\eta_i}:
B'\otimes_A\Frac\,A\to B\otimes_A\Frac\,A.
$$
On the other hand, by \cite[Ch.IV, Prop.17.5.7]{EGA4},
$B$ (resp. $B'$) is the normalization of $A$ in
$B\otimes_A\Frac\,A$ (resp. in $B'\otimes_A\Frac\,A$).
It follows that $h_\eta^\natural$ restricts to a map
$B'\to B$, and the corresponding morphism $E\to E'$ is
necessarily an extension of $h$. This shows that $j^*$
is fully faithful in this case.

$\bullet$\ \
Next, suppose that assumption (b) holds, and set
$\cE:=B^\sim$, $\cE':=B'{}^\sim$; then $\cE$ and
$\cE'$ are \'etale $\cO_{\!X}$-algebras and locally
free $\cO_{\!X}$-modules of finite type, and $h$
corresponds to a morphism $h^\sharp:\cE\to\cE'$.
Under the current assumptions, $j_*\cO_{\!U}$ is
a quasi-coherent $\cO_{\!X}$-algebra, and
$\cO_{\!X,x}$ is integrally closed in $(j_*\cO_{\!U})_x$, 
for every $x\in X$ (\cite[Prop.8.2.31(i)]{Ga-Ra}).
Likewise, $j_*j^*\cE=\cE\otimes_{\cO_{\!X}}j_*\cO_{\!U}$,
so $\cE'_x$ is the integral closure of $\cO_{\!X,x}$
in $(j_*j^*\cE)_x$, for every $x\in X$. Lastly, $\cE_x$
is integral over $\cO_{\!X,x}$, so the map
$(j_*j^*h^\sharp)_x:(j_*j^*\cE)_x\to(j_*j^*\cE')_x$
restricts to a map $\cE_x\to\cE'_x$, and the assertion
follows.

To proceed, we make the following general remark.

\begin{claim}\label{cl_trivial-Zorni}
Let $Z$ be a scheme, $V_0\subset Z$ an open subset, and
$\phi:E\to V_0$ an object of $\bCov(V_0)$. Suppose that,
for every open subset $V\subset Z$ containing $V_0$, the
pull-back functor $j^*_V:\bCov(V)\to\bCov(V_0)$ is fully
faithful. Let $\cF$ be the family consisting of all the
data $(V,\psi,\alpha)$ where $V\subset X$ is any open subset
with $V_0\subset V$, $\psi:E_V\to V$ is a finite \'etale
morphism, and $\alpha:\psi^{-1}V_0\isom E$ is an isomorphism
of $V_0$-schemes. $\cF$ is preordered (see example
\ref{ex_universe}(iii)) by the relation such that
$(V,\psi,\alpha)\geq(V',\psi',\alpha')$ if and only
if $V'\subset V$ and there is an isomorphism
$\beta:\psi^{-1}V'\isom E_{V'}$ of $V'$-schemes,
such that $\alpha'\circ j^*_{V'}(\beta)=\alpha$.
Then the partially ordered quotient $\cF'$ of $\cF$
admits a supremum (see example \ref{ex_poset-quotient}).
\end{claim}
\begin{pfclaim} It is easily seen that, whereas $\cF$
is not small in our universe $\sU$, the quotient $\cF'$
is {\em essentially small}, {\em i.e.} there is a small
subset of $\cF$ mapping bijectively onto $\cF'$ (verification
left to the reader). Using Zorn's lemma, it is easily seen that
every element of $\cF'$ can be dominated by a maximal element,
and it remains to show that any two maximal elements
$(V,\psi,\alpha)$ and $(V',\alpha',\psi)$ of $\cF$ are
isomorphic; to see this, set $V'':=V\cap V'$ : by assumption,
the isomorphism
$\alpha'{}^{-1}\circ\alpha:\psi^{-1}V_0\isom\psi'{}^{-1}V_0$
extends to an isomorphism of $V''$-schemes
$\psi^{-1}V''\isom\psi'{}^{-1}V''$, using which one can glue $E_V$
and $E_{V'}$ to obtain a datum $(V\cup V',\psi'',\alpha'')$ which is
larger than both our maximal elements, hence it is isomorphic to
both.
\end{pfclaim}

\begin{claim}\label{cl_include-pts}
In the situation of claim \ref{cl_trivial-Zorni}, suppose that
$Z$ is reduced and normal, $V_0$ is dense in $Z$, and either one
of the assumptions (a) or (b) of (ii) hold for $Z$ and the open
immersion $V_0\to Z$.
Let $(V_{\max},\psi,\alpha)$ be a supremum for $\cF$, and $\bar z$
a geometric point of $Z$, such that $\phi\times_ZZ(\bar z)$ 
extends to a finite \'etale covering of $Z(\bar z)$; then the support
of $\bar z$ lies in $V_{\max}$.
\end{claim}
\begin{pfclaim} By lemma \ref{lem_go-to-pseudo-lim}, there exist
an \'etale neighborhood $g:Y\to Z$ of $\bar z$, with $Y$ affine,
a finite \'etale morphism $\phi_Y:E_Y\to Y$, and an isomorphism
$h:\phi\times_{V_0}Y\simeq\phi_Y\times_ZV_0$.
We have a natural essentially commutative diagram :
$$
\xymatrix{
\bCov(gY) \ar[r]^-\alpha \ar[d] & \Desc(\bCov,g) \ar[r] \ar[d] &
\bCov(Y) \ar[d] \\
\bCov(V_0\cap gY) \ar[r]^-\beta & \Desc(\bCov,g\times_ZV_0) \ar[r] &
\bCov(Y\times_ZV_0)
}$$
where $\alpha$ and $\beta$ are equivalences, by lemma
\ref{lem_fibred-cat-cov}. Moreover, let $Y':=Y\times_ZY$ and
$Y'':=Y'\times_ZY$; it is easily seen that is assumption
(a) (resp. (b)) of (ii) holds for $Z$ and $V_0$, the same
holds for $Y'$ and $Y'\times_ZV_0$, and also for
$Y''$ and $Y''\times_ZV_0$. By the foregoing, it follows
that the functors
$$
\bCov(Y')\to\bCov(Y'\times_ZV_0)
\qquad\text{and}\qquad
\bCov(Y'')\to\bCov(Y''\times_ZV_0)
$$
are fully faithful, hence the right square subdiagram is
$2$-cartesian (corollary \ref{cor_first-theor}(iii)). Thus, the
datum $(\phi\times_ZgY,\phi_Y,h)$ determines an object $\phi'$
of $\bCov(gY)$, together with an isomorphism
$\phi'\times_ZV_0\simeq\phi\times_ZgY$, which we may use to glue
$\phi$ and $\phi'$ to a single object $\phi''$ of $\bCov(V_0\cup gY)$.
The claim follows.
\end{pfclaim}

The foregoing shows that the assumptions of claim \ref{cl_trivial-Zorni}
are fulfilled, with $Z:=X$, $V_0:=U$ and any object $\phi$
of $\bCov(U)$, hence there exists a largest open subset
$U_{\max}\subset X$ over which $\phi$ extends. However, claim
\ref{cl_include-pts} shows that $U_{\max}=X$, so the proof of
(iii) is complete.

(iv): In view of (iii) and \cite[Ch.IV, Cor.18.8.13]{EGA4},
we are reduced to the case where $X$ is a regular local scheme,
and it suffices to show that $j^*$ is essentially surjective.
We argue by induction on the dimension $n$ of $X\!\setminus\! U$.
If $n=0$, then $X\!\setminus\! U$ is the closed point, in which
case it suffices to invoke the Zariski-Nagata purity theorem
(\cite[Exp.X, Th.3.4(i)]{SGA2}).
Suppose $n>0$ and that the assertion is already known for
smaller dimensions. Let $\phi$ be a given finite \'etale
covering of $U$, and $x$ a maximal point of $X\!\setminus\! U$;
then $X(x)\!\setminus\! U=\{x\}$, so $\phi_{|U\cap X(x)}$ extends
to a finite \'etale morphism $\phi_x$ over $X(x)$.
In turns, $\phi_x$ extends to an affine open neighborhood
$V\subset X$ of $X$, and up to shrinking $V$, this extension
$\phi'$ agrees with $\phi$ on $U\cap V$, by lemma
\ref{lem_go-to-pseudo-lim}. Hence we can glue $\phi$ and $\phi'$,
and replace $U$ by $U\cup V$. Repeating the procedure for every
maximal point of $X\!\setminus\! U$, we reduce the dimension of
$X\!\setminus\! U$; then we conclude by the inductive assumption.
\end{proof}

\begin{definition}\label{def_acycl-loc}
Let $f:X\to S$ be a morphism of schemes, $\bar x$ a
geometric point of $X$ localized at $x\in X$; and set $s:=f(x)$,
$\bar s:=f(\bar x)$. Let also $n\in\N$ be any integer,
and $\L\subset\N$ any non-empty set of prime numbers.
\begin{enumerate}
\item
Denote by $f_{\bar x}:X(\bar x)\to S(\bar s)$ the morphism
of strictly local schemes induced by $f$. We say that $f$
is {\em locally $(-1)$-acyclic at the point $x$}, if the
scheme $f_{\bar x}^{-1}(\xi)$ is non-empty for every strict
geometric point $\xi$ of $S(\bar s)$ (see definition
\ref{def_strict-loc}(i)).
\item
We say that $f$ is {\em locally $0$-acyclic at the point $x$},
if the scheme $f_{\bar x}^{-1}(\xi)$ is non-empty and connected
for every strict geometric point $\xi$ of $S(\bar s)$.
\item
We say that a group $G$ is a {\em finite $\L$-group\/}
if $G$ is finite and all the primes dividing the order
of $G$ lie in $\L$. We say that $G$ is an {\em $\L$-group\/}
if it is a filtered union of finite $\L$-groups.
\item
We say that $f$ is {\em locally $1$-aspherical for $\L$ at
the point $x$}, if we have :
$$
H^1(f_{\bar x}^{-1}(\xi)_\et,G)=\{1\}
$$
for every strict geometric point $\xi$ of $S(\bar s)$, and every
$\L$-group $G$ (where $1$ denotes the trivial $G$-torsor).
\item
We say that $f$ is {\em locally $(-1)$-acyclic} (resp.
{\em locally $0$-acyclic}, resp. {\em locally $1$-aspheri\-cal
for $\L$}), if $f$ is locally $(-1)$-acyclic (resp. locally
$0$-acyclic, resp. locally $1$-aspherical for $\L$) at every
point of $X$.
\item
We say that $f$ is {\em $(-1)$-acyclic} (resp. {\em $0$-acyclic})
if the unit of adjunction : $\cF\to f_*f^*\cF$
is a monomorphism (resp. an isomorphism) for every sheaf
$\cF$ on $S_\et$.
\end{enumerate}
\end{definition}

See section \ref{sec_topoi} for generalities about $G$-torsors for a
group object $G$ on a topos $T$. Here we shall be mainly concerned
with the case where $T$ is the {\'e}tale topos $X^\sim_\et$ of a
scheme $X$, and $G$ is representable by a group scheme, finite and
\'etale over $X$. In this case, using faithfully flat descent one
can show that any $G$-torsor is representable by a {\em principal
$G$-homogeneous space\/}, {\em i.e.} a finite, surjective, {\'e}tale
morphism $E\to X$ with a $G$-action $G\to\Aut_X(E)$ such that the
induced morphism of $X$-schemes
$$
G\times E\to E\times_XE
$$
is an isomorphism.

\sset\subsubsection{}\label{subsec_etale-fungrp}
If $G_X$ is the constant $X^\sim_\et$-group arising from a finite
group $G$ and $X$ is non-empty and connected, right $G_X$-torsors
are also understood as $G$-valued characters of the \'etale
fundamental group of $X$. Indeed, let $\xi$ be a geometric point
of $X$; recall (\cite[Exp.V, \S7]{SGA1}) that $\pi:=\pi_1(X_\et,\xi)$
is defined as the automorphism group of the fibre functor
$$
F_\xi:\bCov(X)\to\Set \qquad (E\xrightarrow{ f } X)\mapsto
f^{-1}\xi.
$$
We endow $\pi$ with its natural profinite topology, as in
\eqref{subsec_Galois-cat-fib}, so that $F_\xi$ can be viewed
as an equivalence of categories
$$
F_\xi:\bCov(X)\to\pi\text{-}\Set.
$$

\begin{lemma}\label{lem_from-coh-to-fundgrp}
In the situation of \eqref{subsec_etale-fungrp}, there exists a
natural bijection of pointed sets :
$$
H^1(X_\et,G_X)\isom H^1_\mathrm{cont}(\pi,G)
$$
from the pointed set of right $G_X$-torsors, to the first non-abelian
continuous cohomology group of $P$ with coefficients in $G$ (see
\eqref{subsec_first-non-abel}).
\end{lemma}
\begin{proof} Let $f:E\to X$ be a right $G_X$-torsor, and fix a
geometric point $s\in F_\xi(E)$; given any $\sigma\in\pi$, there
exists a unique $g_{s,\sigma}\in G$ such that
$$
s\cdot g_{s,\sigma}=\sigma_E(s).
$$
Any $g\in G$ determines a $X$-automorphism $g_E:E\to E$, and by
definition, the automorphism $F_\xi(g_E)$ on $F_\xi(E)$ commutes
with the left action of any element $\tau\in\pi$; however
$F_\xi(g_E)$ is just the right action of $g$ on $F_\xi(E)$, hence we
may compute :
$$
s\cdot g_{s,\tau}\cdot g_{s,\sigma}=(\tau_E(s))\cdot g_{s,\sigma}=
\tau_E(s\cdot g_{s,\tau})=\tau_E\cdot\sigma_E(s)
$$
so the rule $\sigma\mapsto g_{s,\sigma}$ defines a group
homomorphism $\rho_{s,f}:\pi\to G$ which is clearly continuous. We
claim that the conjugacy class of $\rho_{s,f}$ does not depend on
the choice of $s$. Indeed, if $s'\in F_\xi(E)$ is another choice,
there exists a (unique) element $h\in G$ such that $h(s)=s'$;
arguing as in the foregoing we see that $\sigma_E$ commutes with the
right action of $h$ on $F_\xi(E)$. In other words,
$\sigma_E(s)=h^{-1}\circ\sigma_E(s')$, so that $g_{s',\sigma}=h\circ
g_{s,\sigma}\circ h^{-1}$.

Therefore, denote by $\rho_f$ the conjugacy class of $\rho_{s,f}$;
we claim that $\rho_f$ depends only on the isomorphism class of the
$G_X$-torsor $E$. Indeed, any isomorphism $t:E\isom E'$ of right
$G_X$-torsors induces a bijection $F_\xi(t):F_\xi(E)\isom
F_\xi(E')$, equivariant for the action of $G$, and for any
$\sigma\in\pi$ we have $F_\xi(t)\circ\sigma_E=\sigma_{E'}\circ
F_\xi(t)$, whence the assertion.

Conversely, given a continuous group homomorphism $\rho:\pi\to G$,
let us endow $G$ with the induced left $\pi$-action, and right
$G$-action. Then $G$ is an object of $\pi\text{-}\Set$, to which
there corresponds a finite \'etale morphism $E_\rho\to X$, with an
isomorphism $E_\rho\times_X|\xi|\isom G$ of sets with left
$\pi$-action. Since the right action of $G$ is $\pi$-equivariant, we
have a corresponding $G$-action by $X$-automorphisms on $E_\rho$, so
$E_\rho$ is $G$-torsor, and its image under the map of the lemma is
clearly the conjugacy class of $\rho$.

Finally, in order to show that the map of the lemma is injective, it
suffices to prove that, for any right $G_X$-torsor $(f:E\to
X,\phi:E\times G\to E)$ and any $s\in F_\xi(E)$, there exists an
isomorphism of right $G_X$-torsors $E_{\rho_{s,f}}\isom E$. However,
$s$ and $\rho_{s,f}$ determine an identification of sets with left
$\pi$-action :
\set\begin{equation}\label{eq_identify-G-acts}
G\isom F_\xi(E)
\end{equation}
whence an isomorphism $t:E_{\rho_{s,f}}\isom E$ in $\bCov(X)$.
Moreover, \eqref{eq_identify-G-acts} also identifies the right
$G$-action on $F_\xi(E)$ to the natural right $G$-action on $G$; the
latter is $\pi$-equivariant, hence it induces a right $G$-action
$\phi':E\times G\to E$, such that $t$ is $G$-equivariant. To
conclude, it suffices to show that $\phi=\phi'$. In view of
\cite[Ch.IV, Cor.17.4.8]{EGA4}, the latter assertion can be checked
on the stalks over the geometric point $\xi$, where it holds by
construction.
\end{proof}

\begin{remark}\label{rem_from-coh-to-fungrp}
(i)\ \
In the situation of \eqref{subsec_etale-fungrp}, let $E\to X$
be any right $G_X$-torsor, and $\rho_E:\pi\to G$ the corresponding
representation. Then the proof of lemma \ref{lem_from-coh-to-fundgrp}
shows that the left $\pi$-action on $F_\xi(E)$ is isomorphic to the
left $\pi$-action on $G$ induced by $\rho_E$; especially, $\rho_E$
is surjective if and only if $\pi$ acts transitively on $F_\xi(E)$,
if and only if the scheme $E$ is connected (since a decomposition of
$E$ into connected components corresponds to a decomposition of
$F_\xi(E)$ into orbits for the $\pi$-action).

(ii)\ \
Let $\phi:G'\to G$ be a homomorphism of finite groups, and
$$
H^1_\mathrm{cont}(\pi,\phi):
H^1_\mathrm{cont}(\pi,G')\to H^1_\mathrm{cont}(\pi,G)
$$
the induced map. Denote by $r$ (resp. $l$) the right (resp. left)
translation action of $G$ on itself. Let also $E'\to X$ be a
principal $G'$-homogeneous space, given by a map
$\rho:G'\to\Aut_X(E')$, and denote by
$c'\in H^1_\mathrm{cont}(\pi,G')$ the class of $E'$. Then the
class $c:=H^1_\mathrm{cont}(\pi,\phi)(c)$ can be described
geometrically as follows. The scheme $E'\times G$ admits an
obvious right $G$-action, induced by $r$. Moreover, it admits
as well a right $G'$-action : namely, to any element $g\in G'$,
we assign the $X$-automorphism
$\rho_g\times l_{\phi(g^{-1})}:E'\times G\isom E'\times G$.
Set $E:=(E\times G)/G'$; it is easily seen that $E$ is a
principal $G$-homogeneous space, and its class is precisely $c$
(the detailed verification shall be left as an exercise for
the reader).

(iii)\ \
Consider a commutative diagram of schemes
$$
\xymatrix{ E' \ar[rr]^-g \ar[dr]_{f'} & & E \ar[dl]^f \\
           & X
}$$
where $f$ (resp. $f'$) is a $G$-torsor (resp. a $G'$-torsor)
for a given finite group $G$ (resp. $G'$). Denote by
$\rho:G\to\Aut_X(E)$ and $\rho':G'\to\Aut_X(E')$ the
respective actions. We have :
\begin{enumerate}
\alphaenu
\item
Suppose that there exists a group homomorphism
$\phi:G'\to G$, such that $\rho\circ\phi=g\circ\rho'$. Then $g$
induces a $G'$-equivariant map
$$
F_\xi(E')\to\Res(\phi)F_\xi(E)
$$
whence an isomorphism $(E'\times G)/G'\isom G$ of $G$-torsors.
Hence, let $c$ (resp. $c'$) denote any representative of the
equivalence class of $E$ (resp. $E'$) in $H^1_\mathrm{cont}(\pi,G)$
(resp. $H^1_\mathrm{cont}(\pi,G')$); in view of (ii), it follows that
$$
H^1(\pi,\phi)(c')=c.
$$
In other words, the induced diagram of continuous group homomorphisms
$$
\xymatrix{
& \pi_1(X_\et,\xi) \ar[dl]_{c'} \ar[dr]^c \\
G' \ar[rr]^-\phi & & G
}$$
commutes, up to composition with an inner automorphism of $G$.
(Details left to the reader.)
\item
If $E$ is connected, a group homomorphism $\phi:G\to G'$
fulfilling the condition of (a) exists and is unique up
to composition with an inner automorphism of $G$.
Indeed, fix any $e'\in F_\xi(E')$ and let $e:=f(e')$; if $g'\in G'$,
define $\phi(g')$ as the unique $g\in G$ such that
$f(e'\cdot g')=e\cdot g$; also, in view of (i) we may
pick $\sigma_{g'}\in\pi_1(X,\xi)$ such that
$\sigma_{g'}\cdot e'=e'\cdot g'$, and notice that
$\sigma_{g'}\cdot e=e\cdot\phi(g')$ for every $g'\in G'$.
Now, if $h'\in G'$ is any other element, we may compute :
$$
f(e'\cdot g'h')=f(\sigma_{g'}\cdot e'\cdot h')=
\sigma_{g'}\cdot f(e'\cdot h')=\sigma_{g'}\cdot e\cdot\phi(h')=
e\cdot\phi(g')\cdot\phi(h')
$$
whence $\phi(g'h')=\phi(g')\cdot\phi(h')$, as required.
\end{enumerate}
\end{remark}

\begin{lemma}\label{lem_full-faith-0}
Let $f:X\to Y$ be a morphism of schemes, $\cF$, $\cG$ two sheaves on
$Y_\et$. Then :
\begin{enumerate}
\item
If $\cF$ is locally constant and constructible, the natural map :
$$
\theta_f:f^*\cHom_{Y^\sim_\et}(\cF,\cG)\to\cHom_{X^\sim_\et}(f^*\cF,f^*\cG)
$$
is an isomorphism.
\item
If $f$ is $0$-acyclic, it induces a fully faithful functor
$$
f^*:\bCov(Y)\to\bCov(X)\quad :\quad(E\to Y)\mapsto(E\times_YX\to X).
$$
\end{enumerate}
\end{lemma}
\begin{proof}(i): Suppose we have a cartesian diagram of schemes :
$$
\xymatrix{ X' \ar[r]^-{g'} \ar[d]_{f'} & X \ar[d]^f \\
           Y' \ar[r]^-g & Y. }
$$
Then, according to \eqref{eq_obvious-but-cumber}, we have a natural
isomorphism :
$$
\theta_{g'}\circ g^{\prime *}\theta_f \Rightarrow \theta_{f'}\circ
f^{\prime *}\theta_g
$$
(an invertible $2$-cell, in the terminology of
\eqref{sec_2Cats}). Now, if $g$ -- and therefore $g'$ -- is a
covering morphism, $\theta_g$ and $\theta_{g'}$ are isomorphisms,
and $g^{\prime*}\theta_f$ is an isomorphism if and only if the same
holds for $\theta_f$. Summing up, in this case $\theta_f$ is an
isomorphism if and only if the same holds for $\theta_{f'}$. Thus,
we may choose $g$ such that $g^*$ is a constant sheaf, and after
replacing $f$ by $f'$, we may assume that $\cF=S_Y$ is the constant
sheaf associated with a finite set $S$. Since the functors
$$
\cHom_{Y^\sim_\et}(-,\cG):(Y^\sim_\et)^o\to Y^\sim_\et
\quad\text{and}\quad
f^*:Y^\sim_\et\to X^\sim_\et
$$
are left exact, we may further reduce to the case where $S=\{1\}$
is the set with one element, in which case $\cF=1_Y$ is the final
object of $Y^\sim_\et$, and $f^*\cF=1_X$ is the final object of
$X^\sim_\et$. Moreover, we have a natural identification :
$$
\cHom_{Y^\sim_\et}(1_Y,\cG)\isom\cG\quad :\quad\sigma\mapsto\sigma(1)
$$
and likewise for $\cHom_{X^\sim_\et}(1_Y,f^*\cG)$. Using the
foregoing characterization, it is easily checked that, under these
identifications, $\theta_f$ is the identity map of $f^*\cG$, whence
the claim.

(ii): It has already been remarked that $\bCov(Y)$ is equivalent
to the category of locally constant constructible sheaves on $Y_\et$,
and likewise for $\bCov(X)$. Let $E$ and $F$ be two objects of
$\bCov(Y)$; we have natural bijections :
$$
\begin{aligned}
\Hom_{\bCov(X)}(f^*E,f^*F) \isom\, & \Gamma(X,\cHom_{X^\sim_\et}(f^*E,f^*F)) \\
\isom\, & \Gamma(Y,f_*f^*\cHom_{Y^\sim_\et}(E,F)) & & \text{by (i)} \\
\isom\, &
\Gamma(Y,\cHom_{Y^\sim_\et}(E,F)) & & \text{since $f$ is $0$-acyclic} \\
\isom\, & \Hom_{\bCov(Y)}(E,F)
\end{aligned}
$$
as stated.
\end{proof}

\begin{lemma}\label{lem_crit-0-acy}
Let $f:X\to S$ be a quasi-compact morphism of schemes,
and suppose that :
\begin{enumerate}
\alphaenu
\item
$f$ is locally $(-1)$-acyclic.
\item
For every strict geometric point $\xi$ of $S$, the induced morphism
$f_\xi:f^{-1}(\xi)\to|\xi|$ is $0$-acyclic ({\em i.e.} $f$ has
non-empty geometrically connected fibres).
\end{enumerate}
Then $f$ is $0$-acyclic.
\end{lemma}
\begin{proof} Let $\cF$ be a sheaf on $S_\et$. For every
strict geometric point $\xi$ of $S$, we have a commutative diagram :
\set\begin{equation}\label{eq_0-acyclic}
{\diagram
\cF_\xi \ar[r]^-{\eps_\xi} \ar[d]_\alpha & (f_*f^*\cF)_\xi \ar[d] \\
\Gamma(|\xi|,\xi^*\cF) \ar[r]^-{f^*_\xi} &
\Gamma(f^{-1}(\xi),f^*_\xi\circ\xi^*\cF)
\enddiagram}
\end{equation}
where $\eps:\cF\to f_*f^*\cF$ is the unit of adjunction. The map
$\alpha$ is an isomorphism, and the same holds for $f^*_\xi$, since
$f_\xi$ is $0$-acyclic. Hence $\eps_\xi$ is injective, which shows
already that $f$ is $(-1)$-acyclic. It remains to show that
$\eps_\xi$ is surjective. Hence, let $t\in(f_*f^*\cF)_\xi$ be any
section. From \eqref{eq_0-acyclic} we see that there exists a
section $t'\in\cF_\xi$ such that the images of $t$ and
$\eps_\xi(t')$ agree on $\Gamma(f^{-1}(\xi),f^*_\xi\circ\xi^*\cF)$.
We may find an {\'e}tale neighborhood $g:U\to S$ of $\xi$, such that
$t'$ (resp. $t$) extends to a section $t'_U\in\cF(U)$ (resp.
$t_U\in\Gamma(X\times_SU,f^*\cF)$). Let $X_U:=X\times_SU$,
$f_U:=f\times_SU:X_U\to U$, and for every geometric point $\bar x$
of $X_U$, denote by $f^*_{\bar x}:\cF_{f_U(\bar x)}\to
f_U^*\cF_{\bar x}$ the natural isomorphism \eqref{eq_madness}. We
set
$$
V:=\{x\in X_U~|~t_{U,\bar x}=f^*_{\bar x}(t'_{U,f(\bar x)})\}
$$
where $\bar x$ is any geometric point of $X$ localized at $x$, and
$t_{U,\bar x}\in f^*\cF_{\bar x}$ (resp. $t'_{U,f(\bar
x)}\in\cF_{f_U(\bar x)}$) denotes the image of $t_U$ (resp. of
$t'_U$). Clearly $V$ is an open subset of $X_U$, and we have :

\begin{claim}\label{cl_about-V} (i)\ \ $V=f_U^{-1}f_U(V)$.
\begin{enumerate}
\addenu
\item
$f_U(V)\subset U$ is an open subset.
\end{enumerate}
\end{claim}
\begin{pfclaim} (i): Given a point $u\in U$, choose a
strict geometric point $\bar u$ localized at $u$, and set
$\bar s:=g(\bar u)^\rmst$; by assumption, the morphism
$f_{\bar s}:f^{-1}(\bar s)\to|\bar s|$ is $0$-acyclic,
hence the image of $t_U$ in
$\Gamma(f^{-1}(\bar s),f_{\bar s}^*\circ\bar s^*\cF)$
is of the form $f_{\bar s}^*t''$, for some
$t''\in\cF_{\bar s}$.
It follows that $V\cap f^{-1}_U(u)$ is either
the whole of $f^{-1}_U(u)$ or the empty set,
according to whether $t''$ agrees or not with the image
of $t'_U$ in $\cF_{\bar s}=g^*\cF_{\bar u}$.

(ii): The subset $X\!\setminus\! V$ is closed, especially
pro-constructible; since $f$ is quasi-compact, we deduce
that $f_U(X\!\setminus\! V)$ is a pro-constructible
subset of $U$ (corollary \ref{cor-pro-constr}).
It then follows from (i) that $f_U(V)$ is ind-constructible,
hence we are reduced to showing that $f_U(V)$ is closed
under generizations (proposition \ref{prop_closed-under-spec}(ii)).
To this aim, since $V$ is open, it suffices to show
that $f_U$ is {\em generizing\/}, {\em i.e.} that
the induced maps $X_U(x)\to U(u)$ are surjective, for
every $u\in U$ and every $x\in f_U^{-1}(u)$. However,
choose a geometric point $\bar x$ localized at $x$,
and let $\bar u:=f_U(\bar x)$; since the natural maps
$X_U(\bar x)\to X_U(x)$ and $U(\bar u)\to U(u)$ are
surjective, it suffices to show that the same holds
for the map $f_{U,\bar x}:X_U(\bar x)\to U(\bar u)$.
The image of $\bar x$ (resp. $\bar U$) in $X$ (resp.
in $S$) is a geometric point which we denote by the
same name; since the natural maps $X_U(\bar x)\to X(\bar x)$
and $U(\bar u)\to S(\bar u)$ are isomorphisms, we
are reduced to showing that $f_{\bar x}:X(\bar x)\to S(\bar u)$
is surjective, which holds, since $f$ is locally $(-1)$-acyclic.
\end{pfclaim}

Set $W:=f_U(V)$; in view of claim \ref{cl_about-V}, $W$ is an
{\'e}tale neighborhood of $\xi$, and the natural map $\cF(W)\to
f^*\cF(U)$ sends the restriction $t'_{U|W}$ of $t'_U$ to the
restriction $t_{U|V}$ of $t_U$, whence the claim.
\end{proof}

\sset\subsubsection{}\label{subsec_cart-pB}
Consider now a cartesian diagram of schemes :
\set\begin{equation}\label{eq_cart-pB}
{\diagram X' \ar[r]^-{g'} \ar[d]_{f'} & X \ar[d]^f \\
           S' \ar[r]^-g & S
\enddiagram}
\end{equation}
where $g$ is a local morphism of strictly local schemes, and denote
by $s$ (resp. by $s'$) the closed point of $S$ (resp. of $S'$). Let
$x'\in f^{\prime -1}(s')$ be any point, $\bar x{}'$ a geometric
point of $X'$ localized at $x'$, and set $x:=g'(x')$, $\bar
x:=g'(\bar x{}')$. Then $g'$ induces a morphism of $S'$-schemes :
\set\begin{equation}\label{eq_finally}
X'(\bar x{}')\to X(\bar x)\times_SS'.
\end{equation}

\begin{lemma}\label{lem_int-cart-pB}
In the situation of \eqref{subsec_cart-pB}, suppose
that $g$ is an integral morphism. Then :
\begin{enumerate}
\item
The induced morphism $f^{\prime -1}(s')\to f^{-1}(s)$
induces a homeomorphism on the underlying topological
spaces.
\item
\eqref{eq_finally} is an isomorphism.
\end{enumerate}
\end{lemma}
\begin{proof} If $g$ is integral, $\kappa(s')$ is a purely
inseparable algebraic extension of $\kappa(s)$, hence the
morphism $T':=\Spec\,\kappa(s')\to T:=\Spec\,\kappa(s)$
is radicial, and the same holds for the induced
morphisms :
$$
f^{\prime -1}(s')\isom f^{-1}(s)\times_TT'\to f^{-1}(s)
\qquad
f^{-1}_x(s)\times_TT'\to f^{-1}_x(s)
$$
(\cite[Ch.I, Prop.3.5.7(ii)]{EGAI}). Especially (i) holds,
and therefore the natural map $X'(x')\to X(x)\times_SS'$ is
an isomorphism; we see as well that $f^{-1}_x(s)\times_TT'$
is a local scheme. Then the assertion follows from
\cite[Ch.IV, Rem.18.8.11]{EGA4}.
\end{proof}

\begin{proposition}\label{prop_q-fin=basech}
Let $f:X\to S$ and $g:S'\to S$ be morphisms of schemes, with $g$
quasi-finite, and set $X':=X\times_SS'$. Suppose that $f$ is locally
$(-1)$-acyclic (resp. locally $0$-acyclic); then the same holds for
$f':=f\times_SS':X'\to S'$.
\end{proposition}
\begin{proof} Let $\bar s{}'$ be any geometric point of $S'$,
and set $\bar s:=g(\bar s{}')$. Denote by $s\in S$ (resp. $s'\in
S'$) the support of $\bar s$ (resp. $\bar s{}'$); then $f$ is
locally $(-1)$-acyclic at the points of $f^{-1}(s)$, if and only if
$f_{\bar s}:=f\times_SS(\bar s):X\times_SS(\bar s)\to S(\bar s)$
enjoys the same property at the points of $f^{-1}_{\bar s}(\bar s)$.
Likewise, $f'$ is locally $(-1)$-acyclic (resp. locally $0$-acyclic)
at the points of $f^{\prime -1}(s')$, if and only if $f'_{\bar
s{}'}:X'\times_{S'}S'(\bar s{}')\to S'(\bar s{}')$ enjoys the same
property at the points of $f^{-1}_{\bar s}(\bar s)$. Hence, we may
replace $g$ by $g_{\bar s{}'}:S'(\bar s{}')\to S(\bar s)$, and $f$
by the induced morphism $X\times_SS(\bar s)\to S(\bar s)$, which
allows to assume that $g$ is finite (\cite[Ch.IV,
Th.18.5.11]{EGA4}), hence integral. Let $\xi'$ (resp. $\bar x{}'$)
be any strict geometric point of $S'$ (resp. of $f^{\prime
-1}(s')$), and let $\xi:=g(\xi')^\rmst$, (resp. let $\bar x$ be the
image of $\bar x{}'$ in $X$); we have natural morphisms :
$$
f^{\prime -1}_{\bar x{}'}(\xi')\xrightarrow{\alpha}
X(\bar x)\times_S\xi'\xrightarrow{\beta} f^{-1}_{\bar x}(\xi).
$$
However, $\alpha$ is an isomorphism, by lemma
\ref{lem_int-cart-pB}(ii), and $\beta$ is a radicial
morphism, since the field extension $\kappa(\xi)\subset\kappa(\xi')$
is purely inseparable (\cite[Ch.I, Prop.3.5.7(ii)]{EGAI}).
The claim follows.
\end{proof}

\begin{lemma}\label{lem_long-lem}
{\em(i)}\ \
Let $S$ be a strictly local scheme, $s\!\in\!S$ the closed
point, $f\!:\!X\!\to\!S$ a morphism of schemes, $\bar x$
(resp. $\xi$) a strict geometric point of $f^{-1}(s)$
(resp. of $S$). We may find :
\begin{enumerate}
\alphaenu
\item
A cartesian diagram \eqref{eq_cart-pB}, with $S'$ strictly
local, irreducible and normal.
\item
A strict geometric point $\bar x'$ of $f^{\prime -1}(s')$ with
$g'(\bar x{}')^\rmst=\bar x$.
\item
A strict geometric point $\xi'$ of $S'$ localized at the generic
point of $S'$, with $g(\xi')^\rmst=\xi$, and such that
\eqref{eq_finally} induces an isomorphism :
\set\begin{equation}\label{eq_next-finally}
f^{\prime -1}_{\bar x{}'}(\xi')\isom f^{-1}_{\bar x}(\xi).
\end{equation}
\end{enumerate}

{\em(ii)}\ \
Moreover, if $S$ is noetherian, we may find $S'$ as in {\em(i)},
such that $\cO_S(S)$ is a Krull domain.
\end{lemma}
\begin{proof}(i): Denote by $Z\subset S$ the closure of
the image of $\xi$, endow $Z$ with its reduced subscheme
structure, set $Y:=X\times_SZ\subset X$, and let
$h_{\bar x}:Y(\bar x)\to Z$ be the natural morphism. Then $Z$ is a
strictly local scheme (\cite[Ch.IV, Prop.18.5.6(i)]{EGA4}).
Moreover, the closed immersion $Y\to X$ induces an
isomorphism of $Z$-schemes : $Y(\bar x)\isom X(\bar x)\times_SZ$
(lemma \ref{lem_int-cart-pB}(ii)). By construction, $\xi$
factors through a strict geometric point $\xi'$ of $Z$, and
we deduce an isomorphism :
$h_{\bar x}^{-1}(\xi')\isom f^{-1}_{\bar x}(\xi)$
of $Z$-schemes. Thus, we may replace $(S,X,\xi)$ by
$(Z,Y,\xi')$, and assume that $S$ is the spectrum of a
strictly local domain, and $\xi$ is localized at the
generic point of $S$. Say that $S=\Spec\,A$, and denote by
$A^\nu$ the normalization of $A$ in its field of fractions
$F$. Then $A^\nu$ is the union of a filtered family
$(A_\lambda~|~\lambda\in\Lambda)$ of finite $A$-subalgebras
of $F$; since $A$ is henselian, each $A_\lambda$ is a
product of henselian local rings, hence it is a local
henselian ring, so the same holds for $A^\nu$. Moreover,
the residue field $\kappa(s')$ of $A^\nu$ is an algebraic
extension of the residue field $\kappa(s)$ of $A$, which is
separably closed, hence $\kappa(s')$ is separably closed,
{\em i.e.} $A^\nu$ is strictly henselian, so we may fulfill
condition (a) by taking $S':=\Spec\,A^\nu$.
Condition (b) holds as well, due to lemma \ref{lem_int-cart-pB}(i).
Finally, it is clear that $\xi$ lifts to a unique strict
geometric point $\xi'$ of $S'$, and it follows from lemma
\ref{lem_int-cart-pB}(ii) that \eqref{eq_next-finally} is an
isomorphism, as required.

(ii): A direct inspection of the proof of (i) reveals that
if $S$ is noetherian, the scheme $S'$ exhibited is the
spectrum of the normalization of a noetherian domain; the
assertion then follows from \cite[Th.33.10]{Na}.
\end{proof}

\sset\subsubsection{}\label{subsec_lim-strict}
In the situation of \eqref{subsec_filtered-co-cov},
let $x\in X$ be any point, and $s:=f(x)$. Let also $\xi$ be
a strict geometric point of $S$. We deduce a compatible system
of points $x_\lambda:=p'_\lambda(x)\in X_\lambda$, whence a
cofiltered system of local schemes
$$
\cX:=(X_\lambda(x_\lambda)~|~\lambda\in\Lambda).
$$
Moreover, we get a compatible system of strict geometric points
$(\xi_\lambda:=p_\lambda(\xi)^\rmst~|~\lambda\in\Lambda)$, with :
$$
\xi\isom\lim_{\lambda\in\Lambda}\xi_\lambda.
$$
Choose a geometric point $\bar x$ of $X$ localized
at $x$, and set likewise $\bar x_\lambda:=p'_\lambda(\bar x)$;
then $\cX$ lifts to a system
$\cX^\sh:=(X_\lambda(\bar x_\lambda)~|~\lambda\in\Lambda)$,
whose limit is naturally isomorphic to $X(\bar x)$
(\cite[Ch.IV, Prop.18.8.18(ii)]{EGA4}). Furthermore,
$\cX^\sh$ induces a natural isomorphism of $\kappa(\xi)$-schemes :
\set\begin{equation}\label{eq_lim-strict}
f_{\bar x}^{-1}(\xi)\isom
\lim_{\lambda\in\Lambda}f_{\lambda,\bar x_\lambda}^{-1}(\xi_\lambda)
\end{equation}
where, as usual, $f_{\bar x}:X(\bar x)\to S$ (resp.
$f_{\lambda,\bar x_\lambda}:X_\lambda(\bar x_\lambda)\to S_\lambda$)
is deduced from $f$ (resp. from $f_\lambda$). These remarks,
together with the following lemma \ref{lem_excell}, and the
previous lemma \ref{lem_long-lem}, will allow in many cases,
to reduce the study of the fibres of $f_{\bar x}$, to the case
where the base $S$ is strictly local, excellent and normal.

\begin{lemma}\label{lem_excell}
Let $S$ be a strictly local normal scheme. Then there exists
a cofiltered family $\cS:=(S_\lambda~|~\lambda\in\Lambda)$
consisting of strictly local normal excellent schemes,
such that :
\begin{enumerate}
\alphaenu
\item
$S$ is isomorphic to the limit of $\cS$.
\item
The natural morphism $S\to S_\lambda$ is dominant for
every $\lambda\in\Lambda$.
\end{enumerate}
\end{lemma}
\begin{proof} Say that $S=\Spec\,A$, and write $A$ as
the union of a filtered family
$\cA:=(A_\lambda~|~\lambda\in\Lambda)$ of excellent noetherian local
subrings, which we may assume to be normal, by \cite[Ch.IV,
(7.8.3)(ii),(vi)]{EGAIV-2}. Proceeding as in
\eqref{subsec_lim-strict}, we choose a compatible family of
geometric points $\bar s_\lambda$ localized at the closed points of
$\Spec\,A_\lambda$, for every $\lambda\in\Lambda$; using these
geometric points, we lift $\cA$ to a filtered family
$(A_\lambda^\sh~|~\lambda\in\Lambda)$ of strict henselizations,
whose colimit is naturally isomorphic to $A$. Moreover, each
$A_\lambda^\sh$ is noetherian, normal and excellent (\cite[Ch.IV,
Prop.18.8.8(iv), Prop.18.8.12(i)]{EGA4} and proposition
\ref{prop_block-buster}(ii)). Let $\eta$ be the generic point of
$S$, $h_\lambda:S\to S_\lambda:=\Spec\,A^\sh_\lambda$ the natural
morphism, and $\eta^\sh_\lambda:=h_\lambda(\eta)$ for every
$\lambda\in\Lambda$. The cofiltered system
$(S_\lambda~|~\lambda\in\Lambda)$ fulfills condition (a). Moreover,
by construction, the image of $\eta^\sh_\lambda$ in
$\Spec\,A_\lambda$ is the generic point $\eta_\lambda$; then
$\eta^\sh_\lambda$ is the generic point of $S_\lambda$, since the
latter is the only point of $S_\lambda$ lying over $\eta_\lambda$.
Hence (b) holds as well.
\end{proof}

\begin{proposition}\label{prop_acycl}
Let $f:X\to S$ be a flat morphism of schemes. We have :
\begin{enumerate}
\item
$f$ is locally $(-1)$-acyclic.
\item
Suppose moreover, that $f$ has geometrically reduced fibres, and :
\begin{enumerate}
\item
either $f$ is locally finitely presented,
\item
or else, $S$ is locally noetherian.
\end{enumerate}
Then $f$ is locally $0$-acyclic.
\end{enumerate}
\end{proposition}
\begin{proof} Let $x\in X$ be any point, set $s:=f(x)$,
choose a geometric point $\bar x$ of $X$ localized at $x$,
set $\bar s:=f(\bar x)$, and let $\xi$ be any strict
geometric point of $S(\bar s)$.

(i): If $f$ is flat, the induced morphism
$f_{\bar x}:X(\bar x)\to S(\bar s)$ is faithfully flat;
especially, $f_{\bar x}$ is surjective (\cite[Th.7.3(i)]{Mat}).

(ii): Set $X':=X\times_SS(\bar s)$. The natural morphism
$X(\bar x)\to X$ factors uniquely through a morphism of
$S(\bar s)$-schemes $j:X(\bar x)\to X'$, and if we denote
by $\bar x{}'$ the image in $X'$ of $\bar x$, then $j$
induces an isomorphism of
$S(\bar s)$-schemes : $X(\bar x)\isom X'(\bar x{}')$.
Hence, $f$ is locally $0$-acyclic at the point $x$,
if and only if the induced morphism $X'\to S(\bar s)$
is locally $0$-acyclic at the support $x'$ of $\bar x{}'$,
so we may replace $S$ by $S(\bar s)$, and assume that
$S$ is strictly local, when (a) holds, and even
strictly local and noetherian, when (b) holds.
We have to show that $f^{-1}_{\bar x}(\xi)$ is connected,
and by lemma \ref{lem_long-lem}, we are further reduced
to the case where $S=S(\bar s)=\Spec\,A$ is strictly local
and normal, $\xi$ is localized at the generic point of $S$,
and :
\begin{enumerate}
\addenu
\item[(a')]
either $f$ is finitely presented,
\item[(b')]
or else, $A$ is a (not necessarily noetherian) Krull domain.
\end{enumerate}

\begin{claim}\label{cl_covered-above}
In case (b') holds, $f_{\bar x}^{-1}(\xi)$ is connected.
\end{claim}
\begin{pfclaim} Let $F$ be the field of fractions of $A$;
the field $\kappa(\xi)$ is algebraic over $F$, hence
it suffices to show that $X(\bar x)\times_S\Spec\,K$ is
connected for every finite field extension $F\subset K$
(\cite[Ch.IV, Prop.8.4.1(ii)]{EGAIV-3}). Let $A_K$ be
the normalization of $A$ in $K$; then $A_K$ is again
a Krull domain (\cite[Ch.VII, \S1, n.8, Prop.12]{BouAC}),
and $B:=\cO_{\!X,x}^\sh\otimes_AA_K$ is a flat $A_K$-algebra.
Notice that the geometric fibres of the induced morphism
$f_{\bar x,K}:\Spec\,B\to\Spec\,A_K$ are cofiltered limits
of schemes that are {\'e}tale over the fibres of $f$; since
the fibres of $f$ are geometrically reduced, it follows that
the same holds for the fibres of $f_{\bar x,K}$. Hence $B$ is
integrally closed in $B\otimes_AF$ (lemma \ref{lem_Krullu});
especially these two rings have the same idempotents,
whence the contention.
\end{pfclaim}

Finally, suppose that (a') holds. By lemma \ref{lem_excell},
the scheme $S$ is the limit of a cofiltered family
$(S_\lambda~|~\lambda\in\Lambda)$ of strictly local excellent
and normal schemes, such that the natural morphisms
$p_\lambda:S\to S_\lambda$ are dominant. By claim
\ref{cl_descend-geom-reduced}, there exists $\lambda\in\Lambda$
and a flat morphism of schemes $f_\lambda:X_\lambda\to S_\lambda$
with geometrically reduced fibres, with an isomorphism of
$S$-schemes $S\times_{S_\lambda}X_\lambda\isom X$; then, for
every $\mu\in\Lambda$ with $\mu\geq\lambda$, set
$X_\mu:=S_\mu\times_{S_\lambda}X_\lambda$ and
$f_\mu:=S_\mu\times_{S_\lambda}f_\lambda:X_\mu\to S_\mu$. After
replacing $\Lambda$ by a cofinal subset, we may assume
that $f_\mu$ is defined for every $\lambda\in\Lambda$.
For every such $\lambda$, let $x_\lambda\in X_\lambda$ be the
image of $x$. Arguing as in \eqref{subsec_lim-strict}, we obtain
a compatible system of strict geometric points $\xi_\lambda$
of $S_\lambda$ (resp. $\bar x_\lambda$ of $X_\lambda$),
such that $p_\lambda(\xi)$ factors through $\xi_\lambda$;
whence an isomorphism \eqref{eq_lim-strict}.
Thus, $f_{\bar x}^{-1}(\xi)$ is reduced if and only if
$f^{-1}_{\lambda,\bar x_\lambda}(\xi_\lambda)$ is reduced
for every sufficiently large $\lambda\in\Lambda$
(\cite[Ch.IV, Prop.8.7.2]{EGAIV-3}). Furthermore, since
$p_\lambda$ is dominant, $\xi_\lambda$ is localized at the
generic point of $S_\lambda$, for every $\lambda\in\Lambda$.
Thus, we are reduced to the case where $S=S(s)$ is the
spectrum of a strictly local noetherian normal domain
$A$, and $\xi$ is localized at the generic point of $S$;
since such $A$ is a Krull domain (\cite[Th.12.4(i)]{Mat}),
this is covered by claim \ref{cl_covered-above}.
\end{proof}

\begin{example}\label{ex_excel}
(i)\ \
Let $A$ be an excellent local ring, and $A^\wedge$
the completion of $A$. Then the natural morphism :
$$
f:\Spec\,A^\wedge\to\Spec\,A
$$
is locally $0$-acyclic. Indeed, this follows from
proposition \ref{prop_acycl}(ii) (and from the excellence
assumption, which includes the geometric regularity of
the formal fibres of $A$).

(ii)\ \ 
Suppose additionally, that $A$ is strictly local.
Then $f$ is $0$-acyclic. To see this, we apply the criterion
of lemma \ref{lem_crit-0-acy} : indeed, since $f$
is flat, it is $(-1)$-acyclic (proposition \ref{prop_acycl}(i));
it remains to show that $f$ has geometrically connected
fibres, and since $A^\wedge$ is strictly local
(\cite[Ch.IV, Prop.18.5.14]{EGA4}), this is the same as showing
that $f$ is locally $0$-acyclic at the closed point of
$\Spec\,A^\wedge$, which has already been remarked in (i).

(iii)\ \ 
More generally, $f$ is $0$-acyclic whenever
$A$ is excellent and henselian. Indeed, in this case
the argument of (ii) again reduces to showing that
$f$ has geometrically connected fibres. However,
consider the natural commutative diagram :
\set\begin{equation}\label{eq_exc-str-cplte}
{\diagram
\Spec\,(A^\wedge)^\sh \ar[r] \ar[d]_{f^\sh} & \Spec\,A^\wedge \ar[d]^f \\
\Spec\,A^\sh \ar[r] & \Spec\,A.
\enddiagram}
\end{equation}
Since $A$ is henselian, $A^\sh$ is the colimit of a filtered
family of finite {\'e}tale and local $A$-algebras. Since $A$
and $A^\wedge$ have the same residue field, it follows easily
that $A^\wedge\otimes_AA^\sh$ is the colimit of a filtered
family of finite {\'e}tale and local $A^\wedge$-algebras, hence
it is strictly henselian, and therefore \eqref{eq_exc-str-cplte}
is cartesian, especially the geometric fibres of $f$ are
connected if and only if the same holds for the geometric
fibres of $f^\sh$, and the latter are reduced (even regular),
since $A$ is excellent. Hence, we come down to showing that
$f^\sh$ is locally $0$-acyclic at the closed point of
$\Spec\,(A^\wedge)^\sh$, which holds again by proposition
\ref{prop_acycl}(ii).
\end{example}

For future use, we point out the following

\begin{proposition}\label{prop_simpler-conds}
Let $g:X\to Y$ be a flat morphism of excellent noetherian schemes,
with $X$ strictly local, and $Y$ normal. Let $U\subset X$ be
an open subset, and $Z\subset Y$ a closed subscheme. Suppose that :
\begin{enumerate}
\item
$g^{-1}(z)\subset U$ for every maximal point $z$ of $Z$.
\item
$U\cap g^{-1}(z)$ is a dense open subset of $g^{-1}(z)$, for
every $z\in Z$.
\item
The fibres $g^{-1}(z)$ are reduced, for every $z\in Z$.
\end{enumerate}
Then the induced functor $\bCov(U)\to\bCov(U\times_YZ)$ is
fully faithful.
\end{proposition}
\begin{proof} Indeed, say that $X=\Spec\,B$, $Y=\Spec\,A$,
$Z=V(I)$ for some ideal $I\subset A$, and denote by $B^\wedge$
the $\fm_B$-adic completion of the local ring $B$ (where
$\fm_B\subset B$ denotes the maximal ideal).
Let also $f:\Spec\,B^\wedge\to Y$ be the induced morphism, and
$U^\wedge\subset\Spec\,B^\wedge$ the preimage of $U$. In light
of example \ref{ex_excel}(ii) and lemma \ref{lem_full-faith-0}(ii),
it suffices to show that the induced functor
$\bCov(U^\wedge)\to\bCov(U^\wedge\times_YZ)$ is fully faithful.
In view of lemma \ref{lem_Lefschetz-covs}, we are further
reduced to checking that conditions (a)--(c) of proposition
\ref{prop_like-Hartshorne} hold for the induced ring homomorphism
$\phi:A\to B^\wedge$, the open subset $U^\wedge$, and the ideal
$I$. However, by example \ref{ex_Ass}, we have
$\Ass_A(I,A)=\Max(Z)$, hence (c) follows
trivially from our assumption (i). Next, since $B^\wedge$
is a faithfully flat $B$-algebra, assumption (ii) implies that
$U^\wedge\cap f^{-1}(z)$ is a dense open subset of $f^{-1}(z)$,
for every $z\in Z$. Moreover, since $B$ is excellent, the
natural morphism $\Spec\,B^\wedge\to X$ is regular, so the
same holds for the induced morphism $f^{-1}(z)\to g^{-1}(z)$,
and then our assumption (iii) implies -- together with
\cite[Th.32.3(i)]{Mat} -- that $f^{-1}(z)$ is reduced, for every
$z\in Z$, whence condition (b). Lastly, we check condition (a),
{\em i.e.} we show that $B^\wedge$ is $I$-adically complete.
Indeed, let $C$ be the $I$-adic completion
of $B^\wedge$; the natural map $B^\wedge\to C$ is injective,
and it admits a left inverse, constructed as follows.
Let $\underline a:=(a_n~|~n\in\N)$ be a given sequence of elements
of $B^\wedge$, which is Cauchy for the $I$-adic topology; then
$\underline a$ is also Cauchy for the $\fm_B$-adic topology, and
it is easily seen that the limit $l$ of $\underline a$ in the $\fm_B$-adic
topology depends only on the class $[\underline a]$ of $\underline a$
in $C$, so we get a well defined ring homomorphism $\lambda:C\to B^\wedge$
by the rule : $[\underline a]\mapsto l$, and clearly $\lambda$
is the sought left inverse. It remains to check that
$\lambda$ is injective; thus, suppose that $l=0$, and that
$[\underline a]\neq 0$; this means that there exists $N\in\N$
such that $a_n\notin I^N$, for every $n\in\N$. Now, the induced
sequence $(\bar a_n~|~n\in\N)$ of elements of $B^\wedge/I^N$
is stationary, and on the other hand, it converges $\fm_B$-adically
to $0$; therefore $\bar a_n=0$ for every sufficiently large $n\in\N$,
a contradiction.
\end{proof}

\sset\subsubsection{}\label{subsec_comp-formal-coh}
Let $A$ be a noetherian normal ring, and endow the
$A$-algebra $A[[t]]$ with its $t$-adic topology. Let
$$
\varphi:\fX:=\Spf\,A[[t]]\to X:=\Spec\,A[[t]]
\qquad
\pi:X\to S:=\Spec\,A
\qquad
i:S\to X
$$
be respectively the natural morphism of locally ringed spaces,
the natural projection, and the closed immersion determined
by the ring homomorphism $A[[t]]\to A$ given by the rule :
$f(t)\mapsto f(0)$, for every $f(t)\in A[[t]]$. Let also
$U_0\subset\Spec\,A$ be an open subset, $U:=\pi^{-1}U_0$
and $\fU:=\varphi^{-1}U$.
Finally, denote by $\cE$ a locally free $\cO_{\!U}$-module
of finite rank, and set $\cE^\wedge:=\varphi_{|\fU}^*\cE$,
which is a locally free $\cO_\fU$-module of finite rank.

\begin{lemma}\label{lem_ane}
In the situation of \eqref{subsec_comp-formal-coh},
suppose that $S\setminus\! U_0$ has codimension $\geq 2$
in $S$. Then:
\begin{enumerate}
\item
The natural map
$$
\Gamma(U,\cE)\to\Gamma(\fU,\cE^\wedge)
$$
is an isomorphism of $A[[t]]$-modules.
\item
The restriction $i_{|U_0}:U_0\to U$ of\/ $i$ induces an equivalence :
$$
i_{|U_0}^*:\bCov(U)\to\bCov(U_0)
\quad :\quad
(E\to U)\mapsto(E\times_UU_0\to U_0).
$$
\end{enumerate}
\end{lemma}
\begin{proof}(i): To begin with, set $Z:=S\!\setminus\! U_0$;
since the morphism $\pi$ is flat, hence generizing
(\cite[Th.9.5]{Mat}), the closed subset $X\!\setminus\! U=\pi^{-1}Z$
has codimension $\geq 2$ in $X$. Since $A$ and $A[[t]]$
are both normal, we deduce :
\set\begin{equation}\label{eq_serre}
\depth_{X\setminus V}\cO_{\!X}\geq 2
\qquad
\depth_Z\cO_{\!S}\geq 2
\end{equation}
(theorem \ref{th_local-depth} and \cite[Th.23.8]{Mat});
therefore (corollary \ref{cor_local-depth}) :
\set\begin{equation}\label{eq_standard}
\Gamma(U,\cO_{\!X})=\Gamma(X,\cO_{\!X})=A[[t]].
\end{equation}
Next, the short exact sequences of $\cO_{\!X}$-modules :
$$
0\to i_*\cO_{\!S}\to\cO_{\!X}/t^{n+1}\cO_{\!X}\to
\cO_{\!X}/t^n\cO_{\!X}\to 0
\qquad\text{for every $n\in\N$}
$$
induce exact sequences
\set\begin{equation}\label{eq_ane}
R^j\Gamma_{\!\!Z} i_*\cO_{\!S}\to
R^j\Gamma_{\!\!Z}\cO_{\!X}/t^{n+1}\cO_{\!X}\to
R^j\Gamma_{\!\!Z}\cO_{\!X}/t^n\cO_{\!X}
\qquad\text{for every $n,j\in\N$}.
\end{equation}
Then \eqref{eq_serre} and \eqref{eq_ane} yield inductively :
$$
\depth_Z\cO_{\!X}/t^n\cO_{\!X}\geq 2
\qquad\text{for every $n\in\N$}
$$
and again corollary \ref{cor_local-depth} implies :
\set\begin{equation}\label{eq_dept-agin}
\Gamma(U,\cO_{\!X}/t^n\cO_{\!X})=A[t]/t^nA[t]
\qquad\text{for every $n\in\N$}.
\end{equation}
Since $U$ is quasi-compact, we may find a left exact sequence
$P:=(0\to\cE\to\cO_{\!U}^{\oplus m}\to\cO_{\!U}^{\oplus n})$
of $\cO_{\!U}$-modules (corollary \ref{cor_ease-of-ref}).
Since $\varphi$ is a flat morphism of locally ringed spaces,
the sequence $\varphi^*P$ is still left exact.
Since the global section functors are left exact, we
are then reduced to the case where $\cE=\cO_{\!U}$.
Then we may write :
$$
\cE^\wedge=\cO_\fU=\lim_{n\in\N}\,\cO_{\!U}/t^n\cO_{\!U}
$$
where, for each $n\in\N$, we regard $\cO_{\!U}/t^n\cO_{\!U}$
as a sheaf of (pseudo-discrete) rings on $\fU=V(t)\subset U$.
The functor $\Gamma(\fU,-)$ is a right adjoint, hence commutes
with limits, and we deduce an isomorphism :
$$
\Gamma(\fU,\cE^\wedge)\isom
\lim_{n\in\N}\,\Gamma(U,\cO_{\!U}/t^n\cO_{\!U}).
$$
(This is even a homeomorphism, provided we view the target
as a limit of rings with the discrete topology.) Taking
\eqref{eq_dept-agin} into account, we obtain
$\Gamma(\fU,\cE^\wedge)=A[[t]]$ which, together
with \eqref{eq_standard}, implies the contention.

(ii): Notice that (i) and lemma \ref{lem_Lef} imply
that $\Lef(U,i(U_0))$ holds (see definition \ref{def_Lef}).
Since the pull-back functor $\pi_{|U}^*:\bCov(U_0)\to\bCov(U)$
is a right quasi-inverse to $i^*_{|U_0}$, the latter is
essentially surjective. The full faithfulness is a special
case of lemma \ref{lem_Lefschetz-covs}.
\end{proof}

\subsection{Local asphericity of smooth morphisms of schemes}
\label{sec_bcov}
Let $S$ be a strictly local scheme, $s\in S$ the closed point,
$f:X\to S$ a smooth morphism, $\bar x$ any geometric point of
$f^{-1}(s)$, and denote by $f_{\bar x}:X(\bar x)\to S$ the induced
morphism of strictly local schemes. For any open subset $U\subset S$
we have a base change functor :
\set\begin{equation}\label{eq_bcov}
f_{\bar x}^*:\bCov(U)\to\bCov(f^{-1}_{\bar x}U) \qquad (E\to
U)\mapsto(E\times_Uf_{\bar x}^{-1}U).
\end{equation}

\begin{theorem}\label{th_bcov}
In the situation of \eqref{sec_bcov}, we have :
\begin{enumerate}
\item
The functor \eqref{eq_bcov} is fully faithful.
\item
Suppose moreover that $S$ is excellent and normal, and
that $S\!\setminus\! U$ has codimension $\geq 2$ in $S$. Then
\eqref{eq_bcov} is an equivalence of categories.
\end{enumerate}
\end{theorem}
\begin{proof}(i): In view of lemma \ref{lem_full-faith-0}(ii)
it suffices to show that $f_{\bar x}$ is $0$-acyclic (since
in that case, the same will obviously hold also for its restriction
$f^{-1}_{\bar x}U\to U$). To begin with, $f_{\bar x}$ is locally
$(-1)$-acyclic, by proposition \ref{prop_acycl}(i), hence it
remains only to show that $f$ is locally $0$-acyclic at the
point $x$ (lemma \ref{lem_crit-0-acy}). The latter assertion
follows from proposition \ref{prop_acycl}(ii) and
\cite[Ch.IV, Th.17.5.1]{EGA4}.

(ii): In light of (i), it suffices to show that \eqref{eq_bcov}
is essentially surjective, under the assumptions of (ii).
We argue by induction on the relative dimension
$n$ of $f$. Let $x\in X$ be the support of $\bar x$. We may find
an open neighborhood $U\subset X$ of $x$, and an {\'e}tale morphism
of $S$-schemes $\phi:U\to\A^n_S$ (\cite[Ch.IV, Cor.17.11.4]{EGA4}).
Let $\bar x{}':=\phi(\bar x)$; there
follows an isomorphism of $S$-schemes : $X(\bar x)\isom\A^n_S(\bar x{}')$,
hence we may assume from start that $X=\A^n_S$, and
$f$ is the natural projection. Especially, the theorem holds
for $n=0$. Suppose then, that $n>0$, and that the theorem is
already known when the relative dimension is $<n$. Write
$f=h\circ g$, where
$$
g:X\simeq\A^{n-1}_S\times_S\A^1_S\to\A^1_S
\qquad\text{and}\qquad
h:\A^1_S\to S
$$
are the natural projections; set $\bar x_1:=g(\bar x)$, and
$U_1:=h_{\bar x_1}^{-1}U$, (where
$h_{\bar x_1}:S_1:=\A^1_S(\bar x_1)\to S$ is the morphism induced
by $h$). We have $S_1\!\setminus\! U_1=h_{\bar x_1}^{-1}(S\!\setminus\! U)$,
and since flat maps are generizing (\cite[Th.9.5]{Mat})
we easily see that the codimension of $S_1\!\setminus\! U_1$
in $S_1$ equals the codimension of $S\!\setminus\! U$ in $S$.
From our inductive assumption, we deduce that the base change
functor $\bCov(U_1)\to\bCov(f_{\bar x}^{-1}U)$ is essentially
surjective, and hence it suffices to show that the same holds
for the functor $\bCov(U)\to\bCov(U_1)$. Thus, we are reduced
to the case where $X=\A^1_S$. Suppose now, that
$E\to f_{\bar x}^{-1}U$ is a finite {\'e}tale morphism;
we can write $f_{\bar x}$ as the limit of a
cofiltered family of smooth morphisms
$(f_\lambda:Y_\lambda\to S~|~\lambda\in\Lambda)$,
where each $Y_\lambda$ is an affine {\'e}tale
$\A^1_S$-scheme. Then $f^{-1}_{\bar x}U$ is the limit of the
family $(Y_\lambda\times_SU~|~\lambda\in\Lambda)$.
By \cite[Ch.IV, Th.8.8.2(ii), Th.8.10.5]{EGAIV-3} and
\cite[Ch.IV, Prop.17.7.8]{EGA4}, we may find a $\lambda\in\Lambda$,
a finite {\'e}tale morphism $E_\lambda\to Y_\lambda\times_SU$ and
an isomorphism of $f_{\bar x}^{-1}U$-schemes :
$E_\lambda\times_{Y_\lambda}\A^1_S(\bar x)\isom E$. Denote
by $y\in Y_\lambda$ the image of the closed point of $\A^1_S(\bar x)$,
and by $\bar y$ the geometric point of $Y_\lambda$ obtained as
the image of $\bar x$ (the latter is viewed naturally
as a geometric point of $\A^1_S(\bar x)$); by construction,
$y$ lies in the closed fibre
$Y_0:=Y_\lambda\times_S\Spec\,\kappa(s)$, which is an
{\'e}tale $\A^1_{\kappa(s)}$-scheme, and we may therefore
find a specialization $z\in Y_0$ of $y$, with $z$ a closed
point. Pick a geometric point $\bar z$ of $Y_0$ localized at
$z$, and a strict specialization map
$Y_\lambda(\bar z)\to Y_\lambda(\bar y)$
as in \eqref{subsec_strict-special}; there follows a commutative
diagram :
$$
\xymatrix{
\A^1_S(\bar x)\simeq Y_\lambda(\bar y) \ar[r] \ar[d] &
Y_\lambda(\bar z) \ar[d] \\
Y_\lambda(y) \ar[r] & Y_\lambda(z).
}$$
The finite {\'e}tale covering
$E_\lambda\times_{Y_\lambda}Y_\lambda(y)\to Y_\lambda(y)\times_SU$
lies in the essential image of the functor
$$
\bCov(Y_\lambda(z)\times_SU)\to\bCov(Y_\lambda(y)\times_SU)
\qquad
C\mapsto C\times_{Y_\lambda(z)}Y_\lambda(y).
$$
It follows that $E\to f^{-1}_{\bar x}U$ lies in the essential
image of the functor
$$
\bCov(f_{\lambda,\bar z}^{-1}U)\to\bCov(f_{\bar x}^{-1}U)
\qquad
C\mapsto C\times_{Y(\bar z)}Y(\bar y)\simeq
C\times_{Y(\bar z)}\A^1_S(\bar x)
$$
and therefore it suffices to show that the pull-back functor
$\bCov(U)\to\bCov(f^{-1}_{\bar z}U)$ is essentially surjective.
In other words, we may replace $x$ by $z$, and assume throughout
that $x$ is a closed point of $\A^1_S$.

\begin{claim}\label{cl_exist-T}
Under the current assumptions, we may find a strictly local
normal scheme $T$, with closed point $t$, a finite surjective
morphism $g:T\to S$, and a finite morphism of
$\Spec\,\kappa(s)$-schemes :
$$
\Spec\,\kappa(t)\to\Spec\,\kappa(x).
$$
\end{claim}
\begin{pfclaim} Since $\kappa(x)$ is a finite extension
of $\kappa(s)$, it is generated by finitely many algebraic
elements $u_1,\dots u_n$, and an easy induction allows
to assume that $n=1$. In this case, one constructs first
a scheme $T'$ by taking any lifting of the minimal polynomial of
$u_1$ : for the details, see e.g. \cite[Ch.0, (10.3.1.2.)]{EGAIII},
which shows that the resulting $T$ is local, finite and
flat over $S$, so $T'$ maps surjectively onto $S$. Next,
we may replace $T'$ by its maximal reduced subscheme, which
is still strictly local and finite over $S$. Next, since
$S$ is excellent, the normalization $(T')^\nu$ of $T'$
is finite over $S$ (\cite[Ch.IV, Scholie 7.8.3(vi)]{EGAIV-2});
let $T$ be any irreducible component of $(T')^\nu$; by
\cite[Ch.IV, Prop.18.8.10]{EGA4}, $T$ fulfills all the
sought conditions.
\end{pfclaim}

Choose $g:T\to S$ as in claim \ref{cl_exist-T}; since the
residue field extension $\kappa(s)\to\kappa(t)$ is
algebraic and purely inseparable, there exists a unique
point $x'\in\A^1_S(\bar x)\times_ST$ lying over $t$, and
we may find a unique strict geometric point $\bar x{}'$
of $\A^1_S(\bar x)\times_ST$ localized at $x'$, and
lying over $\bar x$. In view of \cite[Ch.IV, Prop.18.8.10]{EGA4},
there follows a natural isomorphism of $T$-schemes :
$$
\A^1_S(\bar x)\times_ST\isom\A^1_T(\bar x{}').
$$
Denote by
$f_{\bar x{}'}:=f_{\bar x}\times_ST:\A^1_T(\bar x{}')\to T$
the natural projection, and set $U_T:=g^{-1}U$; since the
morphism $g:T\to S$ is generizing (\cite[Th.9.4(ii)]{Mat}),
it is easily seen that $T\setminus\! U_T$ has codimension
$\geq 2$ in $T$.

Let $F:\bCov\to\mathbf{Sch}$ be the fibred category
\eqref{eq_fibred-cat-cov}. We have a natural essentially
commutative diagram of categories :
\set\begin{equation}\label{eq_ess-natura}
{\diagram
\bCov(U) \ar[r] \ar[d] & \Desc(F,g\times_SU) \ar[d]^\delta \\
\bCov(f^{-1}_{\bar x}U) \ar[r] & \Desc(F,g\times_Sf_{\bar x}^{-1}U)
\enddiagram}
\end{equation}
where, for any morphism of schemes $h$, we have denoted by
$\Desc(F,h)$ the category of descent data for the fibred
category $F$, relative to the morphism $h$.

According to lemma \ref{lem_fibred-cat-cov}, the morphism
$g$ is of universal $2$-descent for the fibred
category $F$, so the horizontal arrows in \eqref{eq_ess-natura}
are equivalences. Hence, the theorem will follow, once we know
that $\delta$ is essentially surjective. However, we have :

\begin{claim}\label{cl_I=hate-you}
(i)\ \ Set $U'_T:=U_T\times_ST$ and $U''_T:=U'_T\times_ST$.
The pull-back functors :
$$
\bCov(U'_T)\to\bCov(\A^1_T(\bar x{}')\times_TU'_T)
\qquad
\bCov(U''_T)\to\bCov(\A^1_T(\bar x{}')\times_TU''_T)
$$
are fully faithful.
\begin{enumerate}
\addenu
\item
Suppose that the pull-back functor
$$
\bCov(U_T)\to\bCov(f_{\bar x{}'}^{-1}U_T)
$$
is essentially surjective. Then the same holds for the
functor $\delta$.
\end{enumerate}
\end{claim}
\begin{pfclaim}(i): Let $\bar z{}''$ be any geometric point
of $X'':=\A^1_T\times_ST\times_ST$ whose strict image in $\A^1_T$
is $\bar x{}'$, and let $\bar z{}'$ be the image of $\bar z{}''$
in $X':=\A^1_T\times_ST$; by lemma \ref{lem_int-cart-pB}(ii), the
natural morphisms :
$$
X'(\bar z{}')\to\A^1_T(\bar x{}')\times_ST
\qquad
X''(\bar z{}'')\to\A^1(\bar x{}')\times_ST\times_ST
$$
are isomorphisms (notice that $T\times_ST$ is also strictly
local). Then the claim follows from assertion (i) of
the theorem, applied to the projections $X'\to T\times_ST$
and $X''\to T\times_ST\times_ST$.

(ii): Recall that an object of $\Desc(F,g\times_Sf_{\bar x{}'}^{-1}U)$
consists of a finite {\'e}tale morphism $E'_T\to f^{-1}_{\bar x{}'}U_T$
and a $X'$-isomorphism $\beta':E'_T\times_ST\isom T\times_SE'_T$
fulfilling a cocycle condition on $E\times_ST\times_ST$.
By assumption, $E_T$ descends to a finite {\'e}tale morphism
$E_T\to U_T$; then (i) implies that $\beta'$ descends to
a $U'_T$-isomorphism $\beta:E_T\times_ST\isom T\times_SE_T$,
and the cocycle identity for $\beta'$ descends to a cocycle
identity for $\beta$.
\end{pfclaim}

In view of claim \ref{cl_I=hate-you}, we may replace $(S,U,x)$
by $(T,U_T,x')$, and therefore assume that $x$ is a
$\kappa(s)$-rational point of $\A^1_{\kappa(s)}$. In this case,
any choice of coordinate $t$ on $\A^1_S$ yields a section
$\sigma_{\bar x}:S\to\A^1_S(\bar x)$ of the natural projection,
such that $\sigma_{\bar x}(s)=x$. To conclude the proof of the
theorem, it suffices to show that the pull-back functor :
$$
\bCov(f^{-1}_{\bar x}U)\xrightarrow{\sigma^*_{\bar x}}\bCov(U)
$$
is fully faithful.

Say that $S=\Spec\,A$; then the scheme
$\A^1_S(\bar x)$ is the spectrum of $A\{t\}$, the henselization
of $A[t]$ along the ideal $\fm\{t\}$ generated by $t$ and the
maximal ideal $\fm$ of $A$. Let $A^\wedge$ (resp. $A\{t\}^\wedge$)
be the $\fm$-adic (resp. $\fm\{t\}$-adic) completion of $A$
(resp. of $A\{t\}$), and notice the natural isomorphism:
$$
A\{t\}^\wedge/tA\{t\}^\wedge\isom A^\wedge
$$
(indeed, it is easy to check that $A\{t\}^\wedge\simeq A^\wedge[[t]]$),
whence a natural diagram of schemes :
$$
\xymatrix{
X^\wedge:=\Spec\,A\{t\}^\wedge \ar[r]^-{g'} \ar@<-.5ex>[d]_\pi &
\Spec\,A\{t\} \ar@<-.5ex>[d]_{f_{\bar x}} \\
S^\wedge:=\Spec\,A^\wedge \ar[r]^-g \ar@<-.5ex>[u]_\sigma &
\Spec\,A \ar@<-.5ex>[u]_{\sigma_{\bar x}}
}$$
(where $\pi$ is the natural projection) whose horizontal
arrows commute with both the downward arrows and the upward ones.
Set $U^\wedge:= g^{-1}U$; by example \ref{ex_excel}(ii) and lemma
\ref{lem_full-faith-0}(ii), the pull-back functors
$$
g^*:\bCov(U)\to\bCov(U^\wedge)
\qquad
g^{\prime *}:\bCov(f^{-1}_{\bar x}U)\to\bCov(\pi^{-1}U^\wedge)
$$
are fully faithful. Consequently, we are easily reduced
to showing that the pull-back functor :
$\bCov(\pi^{-1}U^\wedge)\xrightarrow{\sigma^*}\bCov(U^\wedge)$
is an equivalence. The latter holds by lemma \ref{lem_ane}(ii).
\end{proof}

\begin{example}\label{ex_bcov}
As an application of theorem \ref{th_bcov}, suppose that
$K\subset E$ is an extension of separably closed fields,
$V_K$ a geometrically normal and strictly local $K$-scheme,
$U\subset V_K$ an open subset, and $\xi$ a geometric point
of $V_E:=V_K\times_KE$, whose image in $V_K$ is supported
on the closed point. Then the induced functor
$$
\bCov(U)\to\bCov(U\times_{V_K}V_E(\xi))
$$
is fully faithful, and it is an equivalence in case
$V_K\!\setminus\!U$ has codimension $\geq 2$ in $V_K$.

Indeed, let $K^\mathrm{a}$ (resp. $E^\mathrm{a}$) be an
algebraic closure of $K$ (resp. $E$), and choose a homomorphism
$K^\mathrm{a}\to E^\mathrm{a}$ extending the inclusion of $K$
into $E$. Then both $V_{K^\mathrm{a}}:=V_K\times_KK^\mathrm{a}$
and $V_{E^\mathrm{a}}(\xi):=V_E(\xi)\times_EE^\mathrm{a}$ are
still normal and strictly local (lemma \ref{lem_int-cart-pB}(ii)),
and the induced functors
$$
\bCov(U)\to\bCov(U\times_KK^\mathrm{a})
\qquad
\bCov(U\times_{V_K}V_E(\xi))\to
\bCov(U\times_{V_K}V_{E^\mathrm{a}}(\xi))
$$
are equivalences (lemma \ref{lem_replace}(i)). It then
suffices to show that the induced functor
$$
\bCov(U\times_KK^\mathrm{a})\to
\bCov(U\times_{V_K}V_{E^\mathrm{a}}(\xi))
$$
has the asserted properties. Hence, we may replace $K$ by
$K^\mathrm{a}$ and $E$ by $E^\mathrm{a}$, and assume from
start that $K\subset E$ is an extension of algebraically
closed fields. In this case, $E$ can be written as the
colimit of a filtered family $(R_\lambda~|~\lambda\in\Lambda)$
of smooth $K$-algebras; correspondingly, $V_E$ is the limit of
a cofiltered system $(V_\lambda~|~\lambda\in\Lambda)$ of smooth
$V_K$-schemes, and -- by lemma \ref{lem_go-to-pseudo-lim} --
$\bCov(U\times_{V_K}V_E(\xi))$ is the $2$-colimit of the system
of categories
$$
\bCov(U\times_{V_K}V_\lambda(\xi_\lambda))
\qquad
(\lambda\in\Lambda)
$$
(where, for each $\lambda\in\Lambda$, we denote by $\xi_\lambda$
the image of $\xi$ in $V_\lambda$). Now the contention follows
directly from theorem \ref{th_bcov}.
\end{example}

\begin{theorem} Let $f:X\to S$ is a smooth morphism of
schemes, $\L\subset\N$ be a set of primes, and suppose
that all the elements of\/ $\L$ are invertible in $\cO_{\!S}$.
Then $f$ is $1$-aspherical for\/ $\L$.
\end{theorem}
\begin{proof} Let $\bar x$ be any geometric point of $X$,
$\bar s:=f(\bar x)$, and $\bar\eta$ a strict geometric point
of $S(\bar s)$. To ease notation, set $T:=X(\bar x)$,
let $f_{\bar x}:T\to S$ be the natural map, and
$T_{\bar\eta}:=f_{\bar x}^{-1}(\bar\eta)$; we have to show
that $H^1(T_{\bar\eta,\et},G)=\{1\}$ for every $\L$-group $G$.
Arguing as in the proof of proposition \ref{prop_acycl},
we reduce to the case where $S=S(\bar s)$. Then, by lemma
\ref{lem_long-lem}, we can further assume that $S$ is
normal and $\bar\eta$ is localized at the generic point $\eta$
of $S$. By lemma \ref{lem_excell}, $S$ is the limit of
a cofiltered system $(S_\lambda~|~\lambda\in\Lambda)$
of strictly local, normal and excellent schemes, and
as usual, after replacing $\Lambda$ by a cofinal subset,
we may assume that $f$ (resp. $\bar\eta$) descends to a
compatible system of morphisms
$(f_\lambda:X_\lambda\to S_\lambda~|~\lambda\in\Lambda)$,
(resp. of strict geometric points $\bar\eta_\lambda$
localized at the generic point of $S_\lambda$).
By \cite[Ch.IV, Prop.17.7.8(ii)]{EGA4}, there exists
$\lambda\in\Lambda$ such that $f_\mu$ is smooth for every
$\mu\geq\lambda$. Then, in view of
\cite[Exp.VII, Rem.5.14]{SGA4-2} and the isomorphism
\eqref{eq_lim-strict}, we may replace $f$ by $f_\lambda$,
and $\bar\eta$ by $\bar\eta_\lambda$, and assume from start
that $S$ is strictly local, normal and excellent, and $G$
is a finite $\L$-group.

Let $\phi:E_{\bar\eta}\to T_{\bar\eta}$ be a principal
$G$-homogeneous space; we come down to showing that
$E_{\bar\eta}$ has a section $T_{\bar\eta}\to E_{\bar\eta}$.
By \cite[Ch.IV, Prop.17.7.8(ii)]{EGA4} and
\cite[Ch.IV, Th.8.8.2(ii), Th.8.10.5]{EGAIV-3}, we may
find a finite separable extension $\kappa(\eta)\subset L$,
and a principal $G$-homogeneous space
$$
\phi_L:E_L\to T_L:=T\times_S\Spec\,L \qquad
\rho_L:G\to\Aut_{T_L}(E_L)
$$
such that
$$
\phi=\phi_L\times_{\Spec\,L}\Spec\,\kappa(\bar\eta)
\qquad
\rho=\rho_L\times_{\Spec\,L}\Spec\,\kappa(\bar\eta).
$$
Say that $S=\Spec\,A$, denote by $A_L$ the normalization
of $A$ in $L$, and set $S_L:=\Spec\,A_L$. Then $S_L$ is
again normal and excellent (\cite[Ch.IV, (7.8.3)(ii),(vi)]{EGAIV-2}),
and the residue field of $A_L$ is an algebraic extension
of the residue field of $A$, hence it is separably closed,
so $S_L$ is strictly local as well. Thus, we may replace
$S$ by $S_L$, and assume that $E_{\bar\eta}$ descends to
a principal $G$-homogeneous space
$E_\eta\to T_\eta:=f_{\bar x}^{-1}(\eta)$
on $T_\eta$. Next, we may write $\eta$ as the limit of
the filtered system of affine open subsets of $S$, so
that -- by the same arguments -- we find an affine open
subset $U\subset S$ and a principal $G$-homogeneous space
$E_U\to T_U:=f^{-1}_{\bar x}U$, with a $G$-equivariant
isomorphism of $T_\eta$-schemes :
$E_U\times_{T_U}T_\eta\isom E_\eta$. Denote by $D_1,\dots,D_n$
the irreducible components of $S\!\setminus\! U$ which have
codimension one in $S$, and for every $i\leq n$, set
$D'_i:=f_{\bar x}^{-1}D_i$. Let also $\eta_T$ be the
generic point of $T$.

\begin{claim}\label{cl_Abhyankar}
For given $i\leq n$, let $y$ be the generic point of $D_i$,
and $z$ a maximal point of $D'_i$. We have:
\begin{enumerate}
\item
$T$ and $E_U$ are normal schemes, and $T(y)$ is regular.
\item
$D'_i$ is a closed subset of pure codimension one in $T$.
\item
Let $\fm_y$ (resp. $\fm_z$) be the maximal ideal of
$\cO_{\!S,y}$ (resp. of $\cO_{\!T,z}$); then
$\fm_y\cdot\cO_{\!T,z}=\fm_z$.
\item
Let $t\in A$ be any element such that $t\cdot\cO_{\!S,y}=\fm_y$.
Then there exist an integer $m>0$ such that $(m,\chara\,\kappa(s))=1$,
a finite {\'e}tale covering
$$
E_y\to T(y)[t^{1/m}]:=T(y)\times_S\Spec\,A[T]/(T^m-t)
$$
and an isomorphism of $T(y)[t^{1/m}]\times_TT_U$-schemes :
$$
E_y\times_TT_U\isom E_U\times_TT(y)[t^{1/m}].
$$
\end{enumerate}
\end{claim}
\begin{pfclaim} (i): Since $S$ is normal by assumption,
the assertion for $T$ and $E_U$ follows from
\cite[Ch.IV, Prop.17.5.7, Prop.18.8.12(i)]{EGA4}.
Next, set $W:=X(y)$; since $\cO_{\!S,y}$ is a discrete
valuation ring, $W$ is a regular scheme
(\cite[Ch.IV, Prop.17.5.8(iii)]{EGA4}). For any
$w\in T(y)\subset W$, the natural map
$\cO_{T(y),w}\to W(w)^\sh$ is faithfully flat, and
$W(w)^\sh$ is regular (\cite[Ch.IV, Cor.18.8.13]{EGA4}),
therefore $\cO_{T(y),w}$ is regular, by
\cite[Ch.0, Prop.17.3.3(i)]{EGAIV}.

(ii): Say that $\fp\subset A$ is the prime ideal of
height one such that $V(\fp)=D_i$; to ease notation, let
also $B:=\cO^\sh_{\!X,x}$.
Let $\{\fq_1,\dots,\fq_k\}\subset\Spec\,B$ be the
set of maximal points of $D'_i$. Using the fact that flat
morphisms are generizing (\cite[Th.9.5]{Mat}), one verifies
easily that $A\cap\fq_j=\fp$ for every $j\leq k$. Fix
$j\leq k$, and set $\fq:=\fq_j$. Since $A$ is normal,
$A_\fp$ is a discrete valuation ring, hence $\fp A_\fp$
is a principal ideal, say generated by $t\in A_\fp$;
therefore $\fq B_\fq$ is the minimal prime ideal of $B_\fq$
containing $t$, so $\fq B_\fq$ has height at most
one, by Krull's Hauptidealsatz (\cite[Th.13.5]{Mat}).
However, a second application of \cite[Th.9.5]{Mat}
shows that the height of $\fq$ in $B$ cannot be lower
than one, hence $\fq$ has height one, which is the contention.

(iii): From (i) and (ii) we see that $\cO_{\!S,y}$
and $\cO_{T,z}$ are discrete valuation rings; then
the assertion follows easily from \cite[Ch.IV, Th.17.5.1]{EGA4}.

(iv): To begin with, since $T(y)$ is regular, it
decomposes as a disjoint union of connected components,
in natural bijection with the set of maximal points
of $D'_i$. Let $Z\subset T(y)$ be the connected open
subscheme containing $z$; it suffices to show that
there exists a finite {\'e}tale covering :
$$
E_Z\to Z[t^{1/m}]:=Z\times_{T(y)}T(y)[t^{1/m}]
$$
with an isomorphism of $Z[t^{1/m}]\times_TT_U$-schemes :
$E_Z\times_TT_U\isom E_U\times_TZ[t^{1/m}]$.
By (iii) we have $t\cdot\cO_{\!T,z}=\fm_z$.
Notice that $T(z)\times_TT_U=T(\eta_T)$, and
$E_{\eta_T}:=E_U\times_TT(z)$ is a disjoint union of
spectra of finite separable extensions $L_1,\dots,L_k$
of $\kappa(\eta_T)$. Moreover, $E_{\eta_T}$ is a principal
$G$-homogeneous space over $T(\eta_T)$, {\em i.e.} every
$L_j$ is a Galois extension of $\kappa(\eta_T)$, with
Galois group $G_j:=\Gal(L_j/\kappa(\eta_T))\subset G$.
Since $G$ is an $\L$-group, the same holds for $G_j$,
hence $E_U\times_TZ$ is tamely ramified along the divisor
$\bar{\{z\}}$ (the topological closure of $\{z\}\subset Z$),
and the assertion follows from Abhyankar's lemma
\cite[Exp.XIII, Prop.5.2]{SGA1}.
\end{pfclaim}

\begin{claim}\label{cl_there-exis} There exist :
\begin{enumerate}
\alphaenu
\item
a finite dominant morphism $S'\to S$, such that both $S'$
and $T':=T\times_SS'$ are strictly local and normal;
\item
an open subset $U'\subset S'$, such that $S'\!\setminus\! U'$
has codimension $\geq 2$ in $S'$;
\item
a finite {\'e}tale morphism $E'\to T'_{U'}:=T\times_SU'$, with
an isomorphism of $T'_{U'}$-schemes :
$$
E'\times_TT_\eta\simeq E_\eta\times_TT'_{U'}.
$$
\end{enumerate}
\end{claim}
\begin{pfclaim} For every $i\leq n$, let $y_i$ be the
maximal point of $D_i$, and choose $t_i\in A$ whose image
in $\cO_{\!S,y_i}$ generates the maximal ideal. Choose also
$m_i\in\N$ with $(m_i,\chara\,\kappa(s))=1$ and such that
there exists a finite {\'e}tale covering
$E_i\to T(y_i)[t_i^{1/m_i}]$ extending the {\'e}tale
covering $E_U\times_TT(y_i)[t_i^{1/m_i}]$
(claim \ref{cl_Abhyankar}(ii.d)). Let $S'$ be the normalization of
$\Spec\,A[t_1^{1/m_1},\dots,t_s^{1/m_n}]$. Then $S'$ is finite
over $S$, hence it is excellent (\cite[Ch.IV, (7.8.3)(ii,vi)]{EGAIV-2}),
and strictly local (cp. the proof of lemma \ref{lem_long-lem}).
Set $E'_\eta:=E_\eta\times_SS'$, $T':=T\times_SS'$; since the
geometric fibres of $f_{\bar x}$ are connected (proposition
\ref{prop_acycl}(ii)), the same holds for the geometric
fibres of the induced morphism $T'\to S'$, therefore $T'$
is connected, and then it is also strictly local, by the usual
arguments. Notice also that $T'$ is the limit of a cofiltered
family of smooth $S'$-schemes, hence it is reduced and normal
(\cite[Ch.IV, Prop.17.5.7]{EGA4}). Say that
$E'_\eta=\Spec\,C$, $T=\Spec\,B$, $T'=\Spec\,B'$, and let $C'$
be the integral closure of $B'$ in $C$. Notice that
$C\otimes_B\kappa(\eta_T)$ is a finite product of finite separable
extensions of the field $B'\otimes_B\kappa(\eta_T)$, and consequently
the natural morphism $\phi':E_{T'}:=\Spec\,C'\to T'$ is finite
(\cite[\S33, Lemma 1]{Mat}). Define :
$$
U':=\{y\in S'~|~
\phi'\times_{S'}S'(y):E_{T'}(y)\to T'(y)\text{ is {\'e}tale}\}.
$$
Let now $y\in U'$ any point; then $S'(y)$ is the limit
of the cofiltered family $(U_\lambda~|~\lambda\in\Lambda)$
of affine open neighborhoods of $y$ in $S'$, and
$\phi'\times_{S'}S'(y)$ the limit of the system of morphisms
$(\phi'_\lambda:=\phi'\times_{S'}U_\lambda~|~\lambda\in\Lambda)$;
we may then find $\lambda\in\Lambda$ such that $\phi'_\lambda$ is
{\'e}tale (\cite[Ch.IV, Prop.17.7.8(ii)]{EGA4}), hence
$U_\lambda\subset U'$, which shows that $U'$ is open.
Furthermore, from \cite[Ch.IV, Prop.17.5.7]{EGA4}
it follows that $E_U\times_TT'$ is normal, whence an
isomorphism of $T'$-schemes :
$$
E_{T'}\times_SU\simeq E_U\times_TT'\qquad
$$
(cp. the proof of lemma \ref{lem_replace}(iii)) especially,
$U\times_SS'\subset U'$. Likewise, by construction we have
natural morphisms : $T'(y_i)\to T(y_i)[t_i^{1/m_i}]$, and
using the fact that all the schemes in view are normal we
deduce isomorphisms of $T'(y_i)$-schemes :
$$
E_{T'}(y_i)\isom E_i\times_{T(y_i)[t_i^{1/m_i}]}T'(y_i).
$$
Thus, $U'$ contains all the points of $S'$ of codimension
$\leq 1$, since the image in $S$ of any such point lies
in $U\cup\{y_1,\dots,y_n\}$. The morphism
$E':=E_{T'}\times_{S'}U'\to T'_{U'}$ fulfills conditions
(a)-(c).
\end{pfclaim}

Now, choose $S'\to S$, $U'\subset S'$, and $E'\to T'_{U'}$
as in claim \ref{cl_there-exis}; since the corresponding
$T'$ is local, there exists a unique point
$x'\in X':=X\times_SS'$ lying over $x$; pick a geometric
point $\bar x{}'$ of $X'$ localized at $x'$, and lying over $\bar x$;
it then follows from \cite[Ch.IV, Prop.18.8.10]{EGA4} that
the natural morphism $X'(\bar x{}')\to T'$ is an isomorphism.
In such situation, theorem \ref{th_bcov} says that there
exists a finite \'etale covering $E\to U'$ with an
isomorphism of $T'_{U'}$-schemes :
$E\times_{U'}T'_{U'}\isom E'$, whence an isomorphism of
$T_{\bar\eta}$-schemes :
$$
E_{\bar\eta}\simeq
E(\bar\eta)\times_{\Spec\,\kappa(\bar\eta)}T_{\bar\eta}.
$$
Since $\kappa(\bar\eta)$ is separably closed, the
{\'e}tale morphism $E(\bar\eta)\to\Spec\,\kappa(\bar\eta)$
admits a section, hence the same holds for $\phi$, as
claimed.
\end{proof}

\sset\subsubsection{}\label{subsec_relative-purity}
Let $f:X\to S$ be a morphism of schemes, and $j:U\subset X$ an
open immersion such that $U_\eta:=U\cap f^{-1}(\eta)\neq\emptyset$
for every $\eta\in S_{\max}$, where $S_{\max}\subset S$ denotes the
subset of all maximal points of $S$. We deduce a natural
essentially commutative diagram of functors :
$$
\cD(S,f,U)\quad:\quad
{\diagram
\bCov(X) \ar[rr]^-{j^*} \ar[d]_{\prod_\eta\iota_\eta^*} & &
\bCov(U) \ar[d]^{\prod_\eta\iota^*_{\eta|U}} \\
\prod_{\eta\in S_{\max}}\bCov(f^{-1}\eta) \ar[rr]^-{\prod_\eta j_\eta^*}
& & \prod_{\eta\in S_{\max}}\bCov(U_\eta)
\enddiagram}\qquad\qquad\qquad
$$
where $j_\eta:U_\eta\to f^{-1}(\eta)$ is the restriction
of $j$ and $\iota_\eta:f^{-1}(\eta)\to X$ is the natural
immersion. Let us say that $U\subset X$ is {\em fibrewise dense},
if $f^{-1}(s)\cap U$ is dense in $f^{-1}(s)$, for every $s\in S$.
Then we have :

\begin{theorem}\label{th_counterpart}
In the situation of \eqref{subsec_relative-purity}, suppose that
$f$ is smooth, and $U$ is fibrewise dense. The following holds :
\begin{enumerate}
\item
The restriction functor $j^*$ is fully faithful.
\item
The diagram $\cD(S,f,U)$ is $2$-cartesian.
\item
If furthermore, $f^{-1}S_{\max}\subset U$, then $j^*$
is an equivalence.
\end{enumerate}
\end{theorem}
\begin{proof} Assertion (ii) means that the functors $j^*$ and
$\iota_\eta^*$ induce an equivalence $(j,\iota_\bullet)^*$ from
$\bCov(X)$ to the category $\cC(X,U)$ of data
\set\begin{equation}\label{eq_datum}
\underline E:=(\phi,(\psi_\eta,\alpha_\eta~|~\eta\in S_{\max}))
\end{equation}
where $\phi$ (resp. $\psi_\eta$) is an object of $\bCov(U)$
(resp. of $\bCov(f^{-1}\eta)$, for every $\eta\in S_{\max}$),
and $\alpha_\eta:\phi\times_U\Spec\,\kappa(\eta)\isom
\psi_\eta\times_{f^{-1}\eta}U_\eta$ is an isomorphism of
$U$-schemes, for every $\eta\in S_{\max}$ (see example
\ref{ex_2-products}(ii)).
On the basis of this description, it is easily seen
that (i),(ii)$\Rightarrow$(iii). Furthermore, we remark :

\begin{claim}\label{cl_taccagni}
(i)\ \ If $j^*$ is fully faithful, then the same holds for
$(j,\iota_\bullet)^*$.
\begin{enumerate}
\addenu
\item
For every open subset $U'\subset X$ containing $U$, suppose that :
\begin{enumerate}
\item
The pull-back functor $\bCov(U')\to\bCov(U)$ is fully faithful.
\item
If $f^{-1}S_{\max}\subset U'$, the pull-back
functor $\bCov(X)\to\bCov(U')$ is an equivalence.
\end{enumerate}
Then assertion (ii) holds.
\end{enumerate}
\end{claim}
\begin{pfclaim}(i): Since $f^{-1}\eta$ is a normal (even regular)
scheme (\cite[Ch.IV, Prop.17.5.7]{EGA4}), the pull-back functors
$\iota^*_\eta$ are fully faithful (lemma \ref{lem_replace}(iii));
the assertion is an immediate consequence.

(ii): In light of (i), it remains only to check that
$(j,\iota_\bullet)^*$ is essentially surjective.
Thus, let $\phi:E\to U$ be a finite {\'e}tale morphism, such that
$i^*_\eta\phi$ extends to a finite {\'e}tale morphism
$\phi'_\eta:E'_\eta\to f^{-1}(\eta)$, for every maximal
point $\eta\in S$.
By claim \ref{cl_trivial-Zorni}, there is a largest open subset
$U_{\max}$ containing $U$, over which $\phi$ extends to a finite
{\'e}tale morphism $\phi_{\max}$. To conclude, we have to show that
$U_{\max}=X$. However, for any maximal point $\eta$, let
$i_\eta:f^{-1}(\eta)\to X(\eta)$ be the natural closed immersion.
By lemma \ref{lem_replace}(i), $i^*_\eta$ is an equivalence,
hence we may find a finite {\'e}tale morphism
$\phi'_{(\eta)}:E'(\eta)\to X(\eta)$ such that
$i^*_\eta\phi'_{(\eta)}\simeq\phi'_\eta$. By the same
token, we also see that $E'(\eta)\times_{X(\eta)}U(\eta)$
is $U(\eta)$-isomorphic to $E\times_UU(\eta)$.

Next, $S(\eta)$ is the limit of the filtered system $\cV$ of all
open subsets $V\subset S$ with $\eta\in V$, hence lemma
\ref{lem_go-to-pseudo-lim} ensures that we may find $V\in\cV$
and an object $\phi'_V:E'_V\to f^{-1}V$ of $\bCov(f^{-1}V)$
such that $\phi'_V\times_VS(\eta)\simeq\phi'_{(\eta)}$, and
after shrinking $V$, we may also assume (again by lemma
\ref{lem_go-to-pseudo-lim}) that
$E'_V\times_XU$ is $U$-isomorphic to $E\times_SV$. Hence
we may glue $E'$ and $E$ along the common intersection,
to deduce a finite {\'e}tale morphism $E'\to U':=U\cup f^{-1}V$
that extends $\phi$. It follows that
$f^{-1}(\eta)\subset f^{-1}V\subset U_{\max}$. Since $\eta$
is arbitrary, (b) implies that the pull-back functor
$\bCov(X)\to\bCov(U_{\max})$ is an equivalence, especially
$\phi$ lies in the essential image of $j^*$, as claimed.
\end{pfclaim}

\begin{claim}\label{cl_firtcase}
(i)\ \ Suppose that $S$ is noetherian and normal,
and $X$ is separated. Then (i) holds.
\begin{enumerate}
\addenu
\item
If furthermore, $S$ is also excellent, then (ii) holds as well.
\end{enumerate}
\end{claim}
\begin{pfclaim}(i): Under the assumptions of the claim, $X$
is normal and noetherian (\cite[Ch.IV, Prop.17.5.7]{EGA4}),
so (i) follows from lemma \ref{lem_replace}(iii), which also
says -- more generally -- that assumption (a) of claim \ref{cl_taccagni}(ii)
holds in this case, hence in order to show (ii) it suffices to
check that assumption (b) of claim \ref{cl_taccagni}(ii) holds
whenever $U\cup f^{-1}S_{\max}\subset U'\subset X$,
especially $X\!\setminus\! U'$ has codimension $\geq 2$ in $X$.
Suppose first that $S$ is regular; then the same holds for
$X$ (\cite[Ch.IV, Prop.17.5.8]{EGA4}), and the contention
follows from lemma \ref{lem_replace}(iv).

In the general case, let $S_\mathrm{reg}\subset S$ be the
regular locus, which is open since $S$ is excellent, and
contains all the points of codimension $\leq 1$, by Serre's
normality criterion (\cite[Ch.IV, Th.5.8.6]{EGAIV-2}).
Consider the restriction $f^{-1}S_\mathrm{reg}\to S_\mathrm{reg}$
of $f$, and the fibrewise dense open immersion
$j_\mathrm{reg}:
U'\cap f^{-1}S_\mathrm{reg}\subset f^{-1}S_\mathrm{reg}$;
by the foregoing, the functor $j^*_\mathrm{reg}$ is an
equivalence, hence we are easily reduced to showing that
the functor $\bCov(X)\to\bCov(U'\cup f^{-1}S_\mathrm{reg})$
is an equivalence, {\em i.e.} we may assume that
$V:=f^{-1}S_\mathrm{reg}\subset U'$. Moreover, since
the full faithfulness of $j^*$ is already known, we only
need to show that any finite {\'e}tale morphism $\phi:E\to U'$
extends to an object of $\bCov(X)$. To this aim, by lemma
\ref{lem_replace}(iii), it suffices to prove that
$\phi\times_{U'}(X(\bar x)\times_XU')$ extends to an
object of $\bCov(X(\bar x))$, for every geometric point
$\bar x$ of $X$. Let $\bar s:=f(\bar x)$, and denote by
$s\in S$ the support of $\bar s$; by assumption, we may
find a geometric point $\xi$ of $f^{-1}(s)$, whose support
lies in $U'\cap f^{-1}(s)$, and a strict specialization
morphism $X(\xi)\to X(\bar x)$. There follows an essentially
commutative diagram :
$$
\xymatrix{\bCov(X(\bar x)) \ar[r]^-\rho \ar[d]_\delta &
\bCov(X(\bar x)\times_XV) \ar[d]^\gamma &
\bCov(S(\bar s)\times_SS_\mathrm{reg}) \ar[l]_\alpha \ar[ld]^\beta \\
\bCov(X(\xi)) \ar[r]^-\tau & \bCov(X(\xi)\times_XV)
}$$
where $\alpha$ and $\beta$ are both equivalences, by theorem
\ref{th_bcov}(ii); hence $\gamma$ is an equivalence as well.
Moreover, both $\bCov(X(\bar x))$ and $\bCov(X(\xi))$ are
equivalent to the category of finite sets, and $\delta$ is
obviously an equivalence. By construction,
$\gamma(\phi\times_{U'}(X(\xi)\times_XV))$ lies in the essential
image of $\tau$, hence $\phi\times_{U'}(X(\xi)\times_XV)$
lies in the essential image of $\rho$, so say it
is isomorphic to $\rho(\phi')$ for some object $\phi'$
of $\bCov(X(\bar x))$. Using (i) (and \cite[Ch.IV, Prop.18.8.12]{EGA4})
one checks easily that
$\phi'\times_XU'\simeq\phi\times_{U'}(X(\xi)\times_XU')$,
whence the contention.
\end{pfclaim}

\begin{claim}\label{cl_pure-case}
Let $m\in\N$ be any integer. Assertions (i) and (ii) hold if $S$
and $X$ are affine schemes of finite type over $\Spec\,\Z$,
the fibres of $f$ have pure dimension $m$, and furthermore :
\set\begin{equation}\label{eq_end-of-fibre}
\dim f^{-1}(s)\!\setminus\! U<m
\qquad\text{for every $s\in S$}.
\end{equation}
\end{claim}
\begin{pfclaim} Indeed, in this situation, $S$ admits
finitely many maximal points, hence the normalization
morphism $S^\nu\to S$ is integral and surjective. Set :
$$
S_2:=S^\nu\times_SS^\nu
\qquad
U_1:=U\times_SS^\nu\qquad U_2:=U\times_SS_2.
$$
Let $\beta:X_1:=X\times_SS^\nu\to X$, $f_1:X_1\to S^\nu$ and
$j_2:U_2\to X_2:=X\times_SS_2$ be the induced morphisms;
clearly $f_1^{-1}(s')$ has pure dimension $m$ for every
$s'\in S^\nu$, and from \eqref{eq_end-of-fibre} we deduce
that $\dim f_1^{-1}(s')\!\setminus\! U_1<m$, especially, $U_1$
is dense in every fibre of $f_1$; by the same token, $U_2$
is dense in $X_2$.
Then $j_2^*$ is faithful (lemma \ref{lem_replace}(ii)),
and lemma \ref{lem_fibred-cat-cov} and corollary
\ref{cor_first-theor}(ii) imply that $j^*$ is
fully faithful, provided the same holds for the functor
$j_1^*:\bCov(X_1)\to\bCov(U_1)$. In other words, in
order to prove assertion (i), we may replace $(f,U)$
by $(f_1,U_1)$, which allows to assume that $S$ is an
affine normal scheme, and then it suffices to invoke claim
\ref{cl_firtcase}, to conclude.

Concerning assertion (ii) : by the foregoing, we already
know that $j^*$ is fully faithful, so the same holds
for $(j,\iota_\bullet)^*$ (claim \ref{cl_taccagni}(i)).
To show that $(j,\iota_\bullet)^*$ is essentially surjective,
let $\underline E$ be an object as in \eqref{eq_datum} of the
category $\cC(X,U)$; the normalization morphism induces a
bijection $S^\nu_{\max}\isom S_{\max}:\eta^\nu\mapsto\eta$, and
clearly $\kappa(\eta^\nu)=\kappa(\eta)$ for every $\eta\in S_{\max}$,
whence a datum
$$
\underline E^\nu:=
(\phi_1:=\phi\times_UU_1,(\psi_\eta,\alpha_\eta~|~\eta^\nu\in S^\nu_{\max}))
$$
of the analogous category $\cC(X_1,U_1)$; by claim
\ref{cl_firtcase}(ii), we may find $\phi'_1\in\Ob(\bCov(X_1))$
and an isomorphism $\alpha:\phi'_1\times_{X_1}U_1\isom\phi_1$. Let
$U_3:=U_2\times_UU_1$, and denote by $j_3:U_3\to X_3:=X_2\times_XX_1$
the natural open immersion; by the foregoing, we know already
that both $j^*_2$ and $j^*_3$ are fully faithful; then corollary
\ref{cor_first-theor}(iii) says that the natural essentially
commutative diagram :
$$
\xymatrix{
\Desc(\bCov,\beta) \ar[r] \ar[d] &
\Desc(\bCov,\beta\times_XU) \ar[d] \\
\bCov(X_1) \ar[r] & \bCov(U_1)
}$$
is $2$-cartesian. Thus, let
$\rho:\bCov(U)\to\Desc(\bCov,\beta\times_XU)$ be the functor
defined in \eqref{subsec_descnt-data}; it follows that
the datum $(\phi'_1,\rho(\phi),\alpha)$ comes from
a descent datum $(\phi'_1,\omega)$ in $\Desc(\bCov,\beta)$.
By lemma \ref{lem_fibred-cat-cov}, the latter descends to
an object $\phi'$ of $\bCov(X)$, and by construction we have
$(j,\iota_\bullet)^*\phi'=\underline E$, as required.
\end{pfclaim}

Next, we consider assertions (i) and (ii) in case where
both $X$ and $S$ are affine. We may find an affine open
covering $X=V_0\cup\cdots\cup V_n$ such that the fibres
of $f_{|V_i}:V_i\to fV_i$ are of pure dimension $i$,
for every $i=0,\dots,n$ (\cite[Ch.IV, Prop.17.10.2]{EGA4}).
For $i=0,\dots,n$, let $j_i:V_i\cap U\to V_i$ be the
induced open immersion; we have natural equivalences
of categories :
$$
\bCov(X)\isom\prod_{i=0}^n\bCov(V_i)
\qquad
\bCov(U)\isom\prod_{i=0}^n\bCov(V_i\cap U)
$$
which induce a natural identification :
$j^*=j^*_0\times\cdots\times j^*_n$.
It follows $j^*$ is fully faithful if and only if the same
holds for every $j_i^*$, and moreover $\cD(S,f,U)$
decomposes as a product of $n$ diagrams $\cD(S,f_{|V_i},U\cap V_i)$.
Hence we may replace $f$ by $f_{|V_m}$, for any $m\leq n$,
after which we may also assume that all the fibres $f$ have
the same pure dimension $m$. In that case, notice that the
assumption on $U$ is equivalent to \eqref{eq_end-of-fibre}.
Next, say that $X=\Spec\,A$, and let $I\subset A$ be an
ideal such that $V(I)=X\!\setminus\! U$; we may write $I$ as
the union of the filtered family $(I_\lambda~|~\lambda\in\Lambda)$
of its finitely generated subideals. Set
$U_\lambda:=X\!\setminus\! V(I_\lambda)$ for every $\lambda\in\Lambda$;
it follows that $U=\bigcup_{\lambda\in\Lambda}U_\lambda$.
For every $\lambda\in\Lambda$, set
$$
Z_\lambda:=\{s\in S~|~\dim f^{-1}(s)\!\setminus\! U_\lambda<m\}.
$$
\begin{claim}\label{cl_little-clumsy}
(i)\ \ $Z_\lambda$ is a constructible subset of $S$,
for every $\lambda\in\Lambda$.
\begin{enumerate}
\addenu
\item
We have $Z_\lambda\subset Z_\mu$ whenever $\mu\geq\lambda$,
and moreover $S=\bigcup_{\lambda\in\Lambda}Z_\lambda$.
\end{enumerate}
\end{claim}
\begin{pfclaim}(i): Let
$V_\lambda:=\{x\in X~|~\dim_x f^{-1}(f(x))\setminus U_\lambda=m\}$;
according to \cite[Ch.IV, Prop.9.9.1]{EGAIV-3}, every $V_\lambda$
is a constructible subset of $X$, hence $f(V_\lambda)$ is a
constructible subset of $S$ (\cite[Ch.IV, Th.1.8.4]{EGAIV}),
so the same holds for $Z_\lambda=S\!\setminus\! f(V_\lambda)$.

(ii): Let $\mu,\lambda\in\Lambda$, such that $\mu\geq\lambda$;
then it is clear that
$f^{-1}(s)\!\setminus\! U_\lambda\subset f^{-1}(s)\!\setminus\! U_\mu$
for every $s\in S$; using \eqref{eq_end-of-fibre}, the claim
follows easily.
\end{pfclaim}

Claim \ref{cl_little-clumsy} and \cite[Ch.IV, Cor.1.9.9]{EGAIV}
imply that $Z_\lambda=S$ for every sufficiently large
$\lambda\in\Lambda$. Hence, after replacing $\Lambda$ by
a cofinal subset, we may assume that all the open subsets
$U_\lambda$ are fibrewise dense. We have a natural essentially
commutative diagram :
$$
\xymatrix{
\bCov(U) \ar[r] \ar[d] &
\Pslim{\lambda\in\Lambda}\,\bCov(U_\lambda) \ar[d] \\
\prod_{\eta\in S_{\max}}\bCov(U_\eta)
\ar[r] & \prod_{\eta\in S_{\max}}
\Pslim{\lambda\in\Lambda}\,\bCov(f^{-1}(\eta)\cap U_\lambda)
}$$
whose horizontal arrows are equivalences
(notation of definition \ref{def_pseudo-lim}(i));
it follows formally that $j^*$ is fully faithful, provided
the same holds for all the pull-back functors
$\bCov(X)\to\bCov(U_\lambda)$, and likewise, $\cD(S,f,U)$ is
$2$-cartesian, provided the same holds for all the
diagrams $\cD(S,f,U_\lambda)$. Hence, we may replace $U$
by $U_\lambda$, and assume that $U$ is constructible, and
\eqref{eq_end-of-fibre} still holds.

Next, we may write $S$ as the limit of a cofiltered family
$(S_\lambda~|~\lambda\in\Lambda)$ of affine schemes of finite
type over $\Spec\,\Z$, and $f$ as the limit of a cofiltered
family $f_\bullet:=(f_\lambda:X_\lambda\to S_\lambda\in\Lambda)$
of affine finitely presented morphisms, such that :
\begin{itemize}
\item
The natural morphism $g_\lambda:S\to S_\lambda$ is dominant
for every $\lambda\in\Lambda$.
\item
$f_\lambda$ is smooth for every $\lambda\in\Lambda$
(\cite[Ch.IV, Prop.17.7.8(ii)]{EGA4}), and
$f_\mu=f_\lambda\times_{S_\lambda}S_\mu$ whenever $\mu\geq\lambda$.
\end{itemize}
Furthermore, we may find $\lambda\in\Lambda$ such that
$U=U_\lambda\times_{S_\lambda}S$ (\cite[Ch.IV, Cor.8.2.11]{EGAIV-3}),
so that $j$ is the limit of the cofiltered system of open
immersions $(j_\mu:U_\mu:=
U_\lambda\times_{S_\lambda}S_\mu\to X_\lambda~|~\mu\geq\lambda)$,
and after replacing $\Lambda$ by a cofinal subset, we may assume
that $j_\mu$ is defined for every $\mu\in\Lambda$.
For every $\lambda\in\Lambda$ and $n\in\N$, let
$X_{\lambda,n}\subset X_\lambda$ be the open and closed subset
consisting of all $x\in X_\lambda$ such that $\dim_xf^{-1}f(x)=n$;
clearly $f_\bullet$ restricts to a cofiltered family
$f_{\bullet,m}:=(f_{\lambda|X_{\lambda,m}}:X_{\lambda,m}\to
S_\lambda~|~\lambda\in\Lambda)$, whose limit is again $f$.
Hence we may replace $X_\lambda$ by $X_{\lambda,m}$, and
assume that the fibres of $f_\lambda$ have pure dimension
$m$, for every $\lambda\in\Lambda$.
For every $\lambda\in\Lambda$, let :
$$
Z'_\lambda:=
\{s\in S_\lambda~|~\dim f_\lambda^{-1}(s)\!\setminus\! U_\lambda=m\}.
$$
and endow $Z'_\lambda$ with its constructible topology $\cT_\lambda$;
since $Z'_\lambda$ is a constructible subset of $S_\lambda$
(\cite[Ch.IV, Prop.9.9.1]{EGAIV-3}), $(Z'_\lambda,\cT_\lambda)$
is a compact topological space, and due to \eqref{eq_end-of-fibre},
we have :
$$
\lim_{\lambda\in\Lambda}Z'_\lambda=\emptyset.
$$
Then \cite[Ch.I, \S9, n.6, Prop.8]{BouTG} implies
that $Z'_\lambda=\emptyset$ for every sufficiently large
$\lambda\in\Lambda$. Set :
$$
\bCov(X_\bullet):=\Pscolim{\mu\geq\lambda}\bCov(X_\mu)
\qquad
\bCov(U_\bullet):=\Pscolim{\mu\geq\lambda}\bCov(U_\mu).
$$
(See definition \ref{def_pseudo-lim}(ii).)
There follows an essentially commutative diagram of categories:
\set\begin{equation}\label{eq_trouble-diag}
{\diagram
\bCov(X) \ar[r] \ar[d]_{j^*} & \bCov(X_\bullet)
\ar[d]^{j^*_\bullet} \\
\bCov(U) \ar[r] & \bCov(U_\bullet)
\enddiagram}
\end{equation}
where $j_\bullet^*$ is the $2$-colimit of the system
of pull-back functors $j_\mu^*:\bCov(X_\mu)\to\bCov(U_\mu)$.
In light of lemma \ref{lem_go-to-pseudo-lim}, the horizontal
arrows of \eqref{eq_trouble-diag} are equivalences, so $j^*$
will be fully faithful, provided the same holds for the functors
$j_\mu^*$, for every large enough $\mu\in\Lambda$.

Hence, in order to prove assertion (i) when $X$ and $S$ are
affine, we may assume that $S$ is of finite type over
$\Spec\,\Z$, the fibres of $f$ have pure dimension $m$, and
\eqref{eq_end-of-fibre} holds, which is the case covered
by claim \ref{cl_pure-case}.

Concerning assertion (ii), since the morphism $g_\mu$ is dominant,
for every $\eta'\in(S_\mu)_{\max}$ we may find $\eta\in S_{\max}$
such that $g_\mu(\eta)=\eta'$. Denote by :
$$
h:f^{-1}\eta\to f^{-1}_\mu\eta'\qquad\text{and}\qquad
j'_\mu:(U_\mu)_{\eta'}:=U_\mu\cap f_\mu^{-1}\eta'\to f^{-1}_\mu\eta'
$$
the natural morphisms. With this notation, we have the following :

\begin{claim}\label{cl_local-eta}
The induced essentially commutative diagram :
$$
\xymatrix{ \bCov(f^{-1}_\mu\eta')
\ar[r]^-{j^{\prime*}_\mu} \ar[d]_{h^*} &
\bCov((U_\mu)_{\eta'}) \ar[d]^{h^*_U} \\
\bCov(f^{-1}\eta) \ar[r]^-{j^*_\eta} & \bCov(U_\eta)
}$$
is $2$-cartesian.
\end{claim}
\begin{pfclaim} The pair $(h^*,j^{\prime*}_\mu)$ induces a
functor $(h,j'_\mu)^*$ from $\bCov(f^{-1}_\mu\eta')$ to the
category of data of the form $(\phi,\phi',\alpha)$, where
$\phi'$ (resp $\phi$) is a finite {\'e}tale covering of
$(U_\mu)_{\eta'}$ (resp. of $f^{-1}\eta$) and
$\alpha:\phi\times_{f^{-1}\eta}U_\eta\isom\phi'\times_{\eta'}\eta$
is an isomorphism in $\bCov(U_\eta)$, and the
contention is that $(h,j'_\mu)^*$ is an equivalence. The full
faithfulness of the functors $j^{\prime*}_\mu$ and $j^*_\eta$
(lemma \ref{lem_replace}(iii)) easily implies the full
faithfulness of $(h,j'_\mu)^*$. To prove that $(h,j'_\mu)^*$
is essentially surjective, amounts to showing that if
$\phi':E'\to(U_\mu)_{\eta'}$ is a finite {\'e}tale morphism and
$$
\phi'':=\phi'\times_{\eta'}\eta:E'':=
E'\times_{\eta'}\eta\to U_\eta
$$
extends to a finite {\'e}tale morphism $\phi:E\to f^{-1}\eta$,
then $\phi'$ extends to a finite {\'e}tale covering of $f^{-1}_\mu\eta'$.
Now, let $L$ be the maximal purely inseparable extension
of $\kappa(\eta')$ contained in $\kappa(\eta)$. Since the
induced morphism $\eta'':=\Spec\,L\to\Spec\,\kappa(\eta')$
is radicial, the base change functors
$$
\bCov(f_\mu^{-1}\eta')\to\bCov((f_\mu^{-1}\eta')\times_{\eta'}\eta'')
\qquad
\bCov((U_\mu)_{\eta'})\to
\bCov((U_\mu)_{\eta'}\times_{\eta'}\eta'')
$$
are equivalences (lemma \ref{lem_replace}(i)).
Thus, we may replace $\eta'$ by $\eta''$, and assume that
the field extension $\kappa(\eta')\subset\kappa(\eta)$ is
separable, hence the induced morphism
$\Spec\,\kappa(\eta)\to\Spec\,\kappa(\eta')$ is regular
(\cite[Ch.VIII, \S7, no.3, Cor.1]{BouAC}), and then the same
holds for the morphism $h$ (\cite[Ch.IV, Prop.6.8.3(iii)]{EGAIV-2}).
Given $\phi'$ as above, set
$\cA:=(j'_{\mu}\circ\phi')_*\cO_{\!E'}$; then $\cA$ is a
quasi-coherent $\cO_{\!f^{-1}_\mu\eta'}$-algebra, and we may
define the quasi-coherent $\cO_{\!f^{-1}_\mu\eta'}$-algebra
$\cB$ as the integral closure of $\cO_{\!f^{-1}_\mu\eta'}$
in $\cA$. By \cite[Ch.IV, Prop.6.14.1]{EGAIV-2}, $h^*\cB$
is the integral closure of $\cO_{\!f^{-1}\eta}$ in
$h^*\cA=j_{\eta*}\circ h^*_U(\phi''_*\cO_{\!E''})=
j_{\eta*}(\phi_*\cO_{\!E})$. By \cite[Ch.IV, Prop.17.5.7]{EGA4},
it then follows that $h^*\cB=\phi^*\cO_{\!E}$,
therefore $\cB$ is a finite {\'e}tale $\cO_{\!f^{-1}_\mu\eta'}$-algebra
(\cite[Ch.IV, Prop.17.7.3(ii)]{EGA4} and
\cite[Ch.IV, Prop.2.7.1]{EGAIV-2}). The claim follows.
\end{pfclaim}

By the foregoing, we already know that $j^*$ is fully faithful,
hence the same holds for $(j,\iota_\bullet)^*$ (claim \ref{cl_taccagni}(i)).
To show that $(j,\iota_\bullet)^*$ is essentially surjective, consider
any $\underline E$ as in \eqref{eq_datum}; we may find $\mu\in\Lambda$,
and a finite {\'e}tale morphism $\phi':E_\mu\to U_\mu$ such that
$\phi=\phi'\times_{U_\mu}U$, whence objects
$\phi'_{\eta'}:=\phi'\times_{U_\mu}(U_\mu)_{\eta'}$
in $\bCov((U_\mu)_{\eta'})$, for every
$\eta'\in(S_\mu)_{\max}$. By construction, we have
$\phi'_{\eta'}\times_{\eta'}\eta\simeq\psi_\eta\times_{f^{-1}\eta}U_\eta$
for every $\eta'\in(S_\mu)_{\max}$ and every $\eta\in S_{\max}$
such that $g_\mu(\eta)=\eta'$. Then claim \ref{cl_local-eta}
shows that, for every $\eta'\in (S_\mu)_{\max}$ there exists
an object $\psi'_{\eta'}$ of $\bCov(f_\mu^{-1}\eta')$ with
isomorphisms :
$$
\psi'_{\eta'}\times_{f^{-1}_\mu\eta'}f^{-1}\eta\simeq\psi_\eta
\qquad
\alpha_{\eta'}:
\psi'_{\eta'}\times_{f^{-1}_\mu\eta'}(U_\mu)_{\eta'}\isom\phi'_{\eta'}.
$$
Therefore, the datum
$\underline E_\mu:=
(\phi',(\psi'_{\eta'},\alpha'_{\eta'}~|~\eta'\in(S_\mu)_{\max}))$
is an object of the $2$-limit of the diagram of categories
$$
\bCov(U_\mu)\xleftarrow{j_\mu^*}\bCov(X_\mu)
\xrightarrow{\prod_{\eta'}\iota_{\eta'}^*}
\prod_{\eta'\in(S_\mu)_{\max}}\bCov(f^{-1}_\mu\eta')
$$
(where $\iota_{\eta'}:f^{-1}_\mu\eta'\to X_\mu$ is the natural
immersion, for every $\eta'\in(S_\mu)_{\max}$). By claim
\ref{cl_pure-case}, the datum $\underline E_\mu$ comes
from an object $\phi_\mu$ of $\bCov(X_\mu)$. Let $\phi''$
be the image of $\phi'_\mu$ in $\bCov(X)$; by construction
we have $(j,\iota_\bullet)^*\phi''=\underline E$, as
required.

This conclude the proof of (i) and (ii), in case $X$ and
$S$ are affine. To deal with the general case, let
$X=\bigcup_{i\in I}V_i$ be a covering consisting of affine
open subschemes, and for every $i\in I$, let
$fV_i=\bigcup_{\lambda\in\Lambda_i}S_{i\lambda}$ be an affine
open covering of the open subscheme $fV_i\subset S$; set also
$V_{i\lambda}:=V_i\cap f^{-1}S_{i\lambda}$ for every $i\in I$
and $\lambda\in\Lambda_i$. The restrictions
$f_{|V_{i\lambda}}:V_{i\lambda}\to S_{i\lambda}$
are smooth morphisms; moreover, the image of the open immersion
$$
j_{i\lambda}:=j_{|U\cap V_{i\lambda}}:U\cap V_{i\lambda}\to V_{i\lambda}
$$
is dense in every fibre of $f_{|V_i}$. The induced morphism
\set\begin{equation}\label{eq_bigcoprod}
\beta:X':=
\coprod_{i\in I}\coprod_{\lambda\in\Lambda_i}V_{i\lambda}\to X
\end{equation}
is faithfully flat, hence of universal $2$-descent for
\eqref{eq_fibred-cat-cov}; moreover, it is easily seen that
$X'':=X'\times_XX'$ is separated, and $j'':=j\times_XX''$ is
a dense open immersion, hence $j''{}^*$ is faithful
(lemma \ref{lem_replace}(ii)). Then, by corollary
\ref{cor_first-theor}(ii), $j^*$ is fully faithful, provided
the pull-back functor $\bCov(X')\to\bCov(X'\times_XU)$ is fully
faithful, {\em i.e.} provided the same holds for the functors
$j^*_{i\lambda}:\bCov(V_{i\lambda})\to\bCov(V_{i\lambda}\cap U)$.
However, each $V_{i\lambda}$ is affine
(\cite[Ch.I, Prop.5.5.10]{EGAI}), hence assertion (i)
is already known  for the morphisms $f_{|V_{i\lambda}}$
and the open subsets $U\cap V_{i\lambda}$; this concludes
the proof of (i).

To show (ii), we use the criterion of claim \ref{cl_taccagni}(ii) :
indeed, assumption (a) is already known, hence we are reduced
to showing that assertion (iii) holds.
To this aim, we consider again the morphism $\beta$ of
\eqref{eq_bigcoprod}, and denote by $f'':X''\to S$ the induced morphism.
Clearly $f''$ is smooth, and $j''$ is an open immersion, such that
$f''{}^{-1}(s)\times_XU$ is dense in $f''^{-1}(s)$, for every
$s\in S$; then assertion (i) implies that $j''{}^*$ is fully
faithful. Moreover it is easily seen that $X''':=X''\times_XX'$
is separated, and $j\times_XX'''$ is a dense open immersion, so
$j'''{}^*$ is faithful (lemma \ref{lem_replace}(ii)), and
therefore corollary \ref{cor_first-theor}(ii) reduces to showing
that the pull-back functor $\bCov(X')\to\bCov(X'\times_XU)$
is an equivalence, or -- what is the same -- that this holds
for the pull-back functors $j^*_{i\lambda}$, which is already
known.
\end{proof}

\sset\subsubsection{}\label{subsec_loc-counterp}
We consider now the local counterpart of theorem \ref{th_counterpart}.
Namely, suppose that $f:X\to S$ is a smooth morphism, let
$\bar x$ be any geometric point of $X$, set $\bar s:=f(\bar x)$,
and let $s\in S$ be the support of $\bar s$.
Then $f$ induces the morphism $f_{\bar x}:X(\bar x)\to S(\bar s)$,
and for every open immersion $j:U\to X(\bar x)$, we may then
consider the diagram $\cD(S(\bar x),f_{\bar x},U)$ as in
\eqref{subsec_relative-purity}. Notice as well that, for every
geometric point $\xi$ of $S$, the fibre $f_{\bar x}^{-1}(\xi)$
is normal, since it is a cofiltered limit of smooth $|\xi|$-schemes;
on the other hand, $f_{\bar x}^{-1}(\xi)$ is also connected,
by proposition \ref{prop_acycl}(ii), hence $f_{\bar x}$ has
geometrically irreducible fibres.

\begin{theorem} In the situation of \eqref{subsec_loc-counterp},
suppose that $U$ contains the generic point of $f^{-1}_{\bar x}(s)$.
Then :
\begin{enumerate}
\item
$j^*:\bCov(X(\bar x))\to\bCov(U)$ is fully faithful.
\item
The diagram $\cD(S(\bar x),f_{\bar x},U)$ is $2$-cartesian.
\end{enumerate}
\end{theorem}
\begin{proof}(i): To begin with, since $f_{\bar x}$ is
generizing (\cite[Th.9.5]{Mat}), and $S$ is local, every
fibre of $f_{\bar x}$ has a point that specializes to the
generic point $\eta_s$ of $f^{-1}_{\bar x}(s)$; since the
fibres are irreducible, it follows that the generic point
of every fibre specializes to $\eta_s$. Therefore $U$ is
fibrewise dense in $X(\bar x)$, and moreover it is connected.
Now, the category $\bCov(X)$ is equivalent to the category
of finite sets, hence every object in the essential image
of $j^*$ is (isomorphic to) a finite disjoint union of copies
of $U$; since $U$ is connected, the morphisms of $U$-schemes
between two such objects $E$ and $E'$ are in natural bijection
with the set-theoretic mappings $\pi_0(E)\to\pi_0(E')$ of
their sets of connected components, whence the assertion.

(ii): Let us write $U$ as the union of a filtered family
$(U_\lambda~|~\lambda\in\Lambda)$ of constructible open
subsets of $X(\bar x)$; up to replacing $\Lambda$ by a
cofinal subset, we may assume that $\eta_s\in U_\lambda$
for every $\lambda\in\Lambda$. Arguing as in the proof
of theorem \ref{th_counterpart}, we see that
$\cD(S(\bar x),f_{\bar x},U)$ is the $2$-limit of
the system of diagrams $\cD(S(\bar x),f_{\bar x},U_\lambda)$,
hence it suffices to show the assertion with $U=U_\lambda$,
for every $\lambda\in\Lambda$, which allows to assume
that $U$ is quasi-compact. Next, arguing as in the proof
of proposition \ref{prop_acycl}, we are reduced to the case
where $S=S(\bar s)$. We may write $X(\bar x)$ as the limit
of a cofiltered system $(X_\lambda~|~\lambda\in\Lambda)$ of
affine schemes, {\'e}tale over $X$, and for $\lambda\in\Lambda$
large enough, we may find an open subset
$U_\lambda\subset X_\lambda$ such that
$U=U_\lambda\times_{X_\lambda}X(\bar x)$
(\cite[Ch.IV, Cor.8.2.11]{EGAIV-3}). For every
$\mu\geq\lambda$, set $U_\mu:=U_\lambda\times_{X_\lambda}X_\mu$,
and denote by $f_\mu:X_\mu\to S$ the natural morphism.
Suppose first that $S$ is irreducible; then, from
lemma \ref{lem_go-to-pseudo-lim} it is easily seen that
$\cD(S,f_{\bar x},U)$ is the $2$-colimit of the
system of diagrams $\cD(S,f_\mu,U_\mu)$, so the assertion
follows from theorem \ref{th_counterpart}(ii) (more generally,
this argument works whenever $S_{\max}$ is a finite set, since
filtered $2$-colimits of categories commute with finite
products).

In the general case, let $\phi:E\to U$ be a finite {\'e}tale
morphism, and suppose that
$\phi_\eta:=\phi\times_{X(\bar x)}f^{-1}_{\bar x}(\eta)$ extends
to an object $\psi_\eta$ of $\bCov(f^{-1}_{\bar x}\eta)$, for
every $\eta\in S_{\max}$. The assertion boils down to showing
that $\phi$ extends to an object $\phi'$ of $\bCov(X(\bar x))$.
To this aim, for every $\eta\in S_{\max}$, let $Z_\eta\to S$
be the closed immersion of the topological closure of $\eta$
in $S$ (which we endow with its reduced scheme structure);
set also $Y_\eta:=X\times_SZ_\eta$. Then $Z_\eta$ is a strictly
local scheme (\cite[Ch.IV, Prop.18.5.6(i)]{EGA4}), and $\bar x$
factors through the closed immersion $Y\to X$, which induces
an isomorphism of $Z$-schemes :
$$
Y_\eta(\bar x)\isom X(\bar x)\times_SZ_\eta
$$
(lemma \ref{lem_int-cart-pB}(ii)). By the foregoing
case, $\phi\times_SZ_\eta$ extends to an object
$\bar\psi_\eta$ of $\bCov(Y_\eta(\bar x))$. However, $Z_\eta$
is the limit of the cofiltered system $(Z_{\eta,i}~|~i\in I(\eta))$
consisting of the constructible closed subschemes of $S$ that
contain $Z_\eta$. By lemma \ref{lem_go-to-pseudo-lim}, it
follows that we may find $i\in I(\eta)$ and an object
$\bar\psi_{\eta,i}$ of $\bCov(X(\bar x)\times_SZ_{\eta,i})$
whose image in $\bCov(Y_\eta(\bar x))$ is isomorphic to
$\bar\psi_\eta$, and if $i$ is large enough,
$\bar\psi_{\eta,i}\times_{X(\bar x)}U$ agrees with
$\phi\times_SZ_{\eta,i}$ in $\bCov(U\times_SZ_{\eta,i})$.
For each $\eta,\eta'\in S_{\max}$, choose $i\in I(\eta)$,
$i'\in I(\eta')$ with these properties, and to ease notation, set :
$$
X'_\eta:=X\times_SZ_{\eta,i}\qquad
X''_{\eta\eta'}:=X'_\eta\times_SZ_{\eta',i'}\qquad
\phi'_\eta:=\bar\psi_{\eta,i}
$$
and denote by
$\alpha_\eta:\phi'_\eta\times_{X(\bar x)}U\isom\phi\times_SZ_{\eta,i}$
the given isomorphism. As in the foregoing,
we notice that $\bar x$ factors through $X'_\eta$, and the
closed immersion $X'_\eta\to X$ induces an isomorphism
$X'_\eta(\bar x)\isom X(\bar x)\times_SZ_{\eta,i}$
of $Z_{\eta,i}$-schemes. According to \cite[Ch.IV, Cor.1.9.9]{EGAIV},
we may then find a finite subset $T\subset S_{\max}$ such that
the induced morphism :
$$
\beta:X_1:=\coprod_{\eta\in T}X'_\eta(\bar x)\to X(\bar x)
$$
is surjective. Set $X_2:=X_1\times_{X(\bar x)}X_1$ and
$X_3:=X_2\times_{X(\bar x)}X_1$; notice that $X_2$ is
the disjoint union of schemes of the form
$X(\bar x)\times_SZ_{\eta,i}\times_SZ_{\eta',i'}$,
for $\eta,\eta'\in T$, and again, the latter is naturally
isomorphic to $X''_{\eta\eta'}(\bar x)$, for a unique lifting
of the geometric point $\bar x$ to a geometric point of $X''_{\eta\eta'}$.
Similar considerations can be repeated for $X_3$, and
in light of (i), we deduce that the pull-back functors :
$$
\bCov(X_i\times_{X(\bar x)}U)\to\bCov(X_i)\qquad i=1,2,3
$$
are fully faithful, in which case corollary \ref{cor_first-theor}(iii)
says that the essentially commutative diagram of categories :
$$
\xymatrix{
\Desc(\bCov,\beta) \ar[r] \ar[d] &
\Desc(\bCov,\beta\times_{X(\bar x)}U) \ar[d] \\
\prod_{\eta\in T}\bCov(X'_\eta(\bar x)) \ar[r] &
\prod_{\eta\in T}\bCov(X'_\eta(\bar x)\times_{X(\bar x)}U)
}$$
is $2$-cartesian.
Let $\rho:\bCov(U)\to\Desc(\bCov,\beta\times_{X(\bar x)}U)$
be the functor defined in \eqref{subsec_descnt-data}; it
follows that the datum
$((\phi'_\eta,\alpha_\eta~|~\eta\in T),\rho(\phi))$ comes from
a descent datum $(\phi'_1,\omega)$ in $\Desc(\bCov,\beta)$.
By lemma \ref{lem_fibred-cat-cov}, the latter descends to
an object $\phi'$ of $\bCov(X(\bar x))$, and by construction
we have $j^*\phi'=\phi$, as required.
\end{proof}

\subsection{\'Etale coverings of log schemes}
We resume the general notation of \eqref{subsec_special-schs},
especially, we choose implicitly $\tau$ to be either the Zariski
or \'etale topology, and we shall omit further mention of this
choice, unless the omission might be a source of ambiguities.

\sset\subsubsection{}\label{subsec_subdiv-and-pi_one} Let
$\underline Y:=((Y,\underline N),T,\psi)$ be an object of
$\cK_\intg$ (see \eqref{subsec_not_quite_fibration} : especially
$\tau=\Zar$ here), $\phi:T'\to T$ an integral proper subdivision
of the fan $T$ (definition \ref{def_subdivisions}(ii),(iii)),
and suppose that both $T$ and $T'$ are locally fine and
saturated. Set
$((Y',\underline N'),T',\psi'):=\phi^*\underline Y$ (proposition
\ref{prop_basis-technique}(iii)), and let $f:Y'\to Y$ be the
morphism of schemes underlying the cartesian morphism
$\phi^*\underline Y\to\underline Y$.

\begin{proposition}\label{prop_from-proper-base-change}
In the situation of \eqref{subsec_subdiv-and-pi_one}, the following
holds :
\begin{enumerate}
\item
For every geometric point $\xi$ of\/ $Y$, the fibre $f^{-1}(\xi)$ is
non-empty and connected.
\item
The functor $f^*:\bCov(Y)\to\bCov(Y')$ is an equivalence.
\end{enumerate}
\end{proposition}
\begin{proof}(i): The assertion is local on $Y$, hence we may assume
that $T=(\Spec\,P)^\sharp$ for a fine, sharp and saturated monoid
$P$. In this case, we may find a subdivision $\phi':T''\to T'$ such
that both $\phi'$ and $\phi\circ\phi'$ are compositions of saturated
blow up of ideals generated by at most two elements of $P$ (example
\ref{ex_cut-by-hyperplane}(iii)). We are reduced to showing the
assertion for the morphisms $\phi'$ and $\phi\circ\phi'$, after
which, we may further assume that $\phi$ is the saturated blow up of
an ideal generated by two elements of $\Gamma(T,\cO_T)$, in which
case the assertion follows from the more precise theorem
\ref{th_elementary-blowup}.

(ii): First, we claim that the assertion is local on the Zariski
topology of $Y$. Indeed, let $Y=\bigcup_{i\in I}U_i$ be a Zariski
open covering, and set $U_{ij}:=U_i\cap U_j$ for every $i,j\in I$;
according to corollary \ref{cor_first-theor}(ii) and lemma
\ref{lem_fibred-cat-cov}, it suffices to prove the contention for
the objects $(U_i\times_Y(Y,\underline N),T,\psi_{|U_i})$ and
$(U_{ij}\times_Y(Y,\underline N),T,\psi_{|U_{ij}})$, and the
morphism $\phi$. Hence, we may again suppose that $T=\Spec\,P$, for
a monoid $P$ as in (i).

Arguing as in the proof of (i), we may next reduce to the case
where $\phi$ is the saturated blow up of an ideal generated by
two elements of $\Gamma(T,\cO_T)$.
Now, the morphism of schemes $f$ induces a morphism of topoi
$Y^{\prime\sim}_\et\to Y^\sim_\et$ which we denote again by $f$;
then, assertion (i) and the proper base change theorem
\cite[Exp.XII, Th.5.1(i)]{SGA4-3} imply that the unit of
adjunction $\cF\to f_*f^*\cF$ is an isomorphism for every
sheaf $\cF$ on $Y$, so $f^*$ is fully faithful (proposition
\ref{prop_fullfaith-adjts}(iii)). Next, if $\phi:E'\to Y'$
is a finite \'etale covering, $Y'$ decomposes into a disjoint
union of open and closed subsets $Y'=Y'_0\cup\cdots\cup Y'_k$,
such that $(\phi_*\cO_{\E})_{|Y'_r}$ is locally free of rank $r$,
for every $r=0,\dots,k$. Set $Y_r:=f(Y'_r)$ for $r=0,\dots,k$;
since the fibers of $f$ are connected, we see that
$Y'_r=f^{-1}Y_r$ for every such $r$. Moreover, since $f$ is a
closed map and the topology of $Y$ is induced from that of
$Y'$ via $f$, we see that $Y_r$ is a closed subset of $Y$
for $r=0,\dots,k$, so $Y=Y_0\cup\cdots\cup Y_k$ is a
partition of $Y$ by open and closed subsets. After replacing
$Y$ by $Y_r$, we may then assume that $\phi_*\cO_{\!E}$ is
locally free of rank $r$, and notice that the isomorphism
classes of \'etale coverings of this type are classified
by the pointed set $H^1(Y'_\et,S_{r,Y'})$, where $S_r$ is the
symmetric group on $r$ elements. Thus, it suffices to show
that the induced map :
\set\begin{equation}\label{eq_down-to-H^1}
H^1(Y_\et,S_{r,Y})\to H^1(Y'_\et,S_{r,Y'})
\end{equation}
is a bijection. However, theorem \ref{th_Leray-H^1} yields
the exact sequence of pointed sets :
$$
\{1\}\to H^1(Y_\et,f_*S_{r,Y'})\xrightarrow{ u } H^1(Y'_\et,S_{r,Y'})
\to H^0(Y_\et,R^1f_*S_{r,Y'}).
$$
On the other hand, assertion (i) and the proper base change
theorem \cite[Exp.XII, Th.5.1(i)]{SGA4-3} imply that the
unit of adjunction $S_{r,Y}\to f_*f^*S_{r,Y}=f_*S_{r,Y'}$ is an
isomorphism, hence $u$ is naturally identified to
\eqref{eq_down-to-H^1}, and we are reduced to showing that
the natural morphism
$$
\tau_{f,S_{r,Y'}}:1_{Y^\sim_\et}\to R^1f_*S_{r,Y'}
$$
is an isomorphism (notation of \eqref{subsec_define-taus}).
The latter can be checked on the stalks, and in view of
\cite[Exp.XII, Cor.5.2(ii)]{SGA4-3} we are reduced to showing
that $H^1(f^{-1}(\xi)_\et,S_r)=\{1\}$ for every geometric point $\xi$
of $Y$. However, according to theorem \ref{th_elementary-blowup},
the reduced geometric fibre $f^{-1}(\xi)_\mathrm{red}$ is either
isomorphic to $|\xi|$ (in which case the contention is trivial),
or else it is isomorphic to the projective line $\P^1_{\kappa(\xi)}$,
in which case -- in view of lemma \ref{lem_replace}(i) -- it suffices
to show that every finite \'etale morphism $E\to\P^1_{\kappa(\xi)}$
admits a section, which is well known.
\end{proof}

The class of \'etale morphisms of log schemes was introduced
in section \ref{sec_smooth-log-sch} : its definition and its
main properties parallel those of the corresponding notion
for schemes, found in \cite[Ch.IV, \S17]{EGAIV}. In the present
section this theme is further advanced : we will consider the
logarithmic analogue of the classical notion of {\em \'etale
covering} of a scheme. To begin with, we make the following :

\begin{definition}\label{def_Kummer-log-sch}
(i)\ \
Let $\phi:T\to S$ be a morphism of fans. We say that $\phi$
is {\em of Kummer type}, if the map
$(\log\phi)_t:\cO_{\!S,\phi(t)}\to\cO_{T,t}$ is of Kummer type,
for every $t\in T$ (see definition \ref{def_Kummer-monoids}).

(ii)\ \
Let $f:(Y,\underline N)\to(X,\underline M)$ be a morphism
of log schemes.
\begin{enumerate}
\alphaenu
\item
We say that $f$ is {\em of Kummer type}, if for every
$\tau$-point $\xi$ of $Y$, the morphism of monoids
$(\log f)_\xi:f^*\underline M_\xi\to\underline N_\xi$
is of Kummer type.
\item
A {\em Kummer chart} for $f$ is the datum of charts
$$
\omega_P:P_Y\to\underline N \qquad \omega_Q:Q_X\to\underline M
$$
and a morphism of monoids $\theta:Q\to P$ such that
$(\omega_P,\omega_Q,\theta)$ is a chart for $f$ (see definition
\ref{def_chart}(iii)), and $\theta$ is of Kummer type.
\end{enumerate}
\end{definition}

\begin{remark}\label{rem_Kummer-type}
(i)\ \
Let $f:P\to Q$ be a morphism of monoids of Kummer type, with $P$
integral and saturated. It follows easily from lemma
\ref{lem_Kummer-fans}(iii,v), that the induced morphism of fans
$(\Spec\,f)^\sharp:(\Spec\,Q)^\sharp\to(\Spec\,P)^\sharp$ is of
Kummer type.

(ii)\ \
In the situation of definition \ref{def_Kummer-log-sch}(ii), it is
easily seen that $f$ is of Kummer type, if and only if the morphism
of $Y$-monoids
$f^*\underline M{}^\sharp_\xi\to\underline N{}^\sharp_\xi$ deduced
from $\log f$, is of Kummer type for every $\tau$-point $\xi$ of $Y$.

(iii)\ \
In the situation of definition \ref{def_Kummer-log-sch}(ii.b),
suppose that the chart $(\omega_P,\omega_Q,\theta)$ is of Kummer
type, and $Q$ is integral and saturated. It then follows easily
from (i), (ii) and example \ref{ex_fan-andlogscheme}(ii), that
$f$ is of Kummer type. In the same vein, we have the following :
\end{remark}

\begin{lemma}\label{lem_needed_once}
Let $f:(Y,\underline N)\to(X,\underline M)$ be a morphism of log
schemes with coherent log structures, $\xi$ a $\tau$-point of\/
$Y$, and suppose that :
\begin{enumerate}
\alphaenu
\item
\ $\underline M_{f(\xi)}$ is fine and saturated.
\item
The morphism $(\log f)_\xi:\underline M_{f(\xi)}\to\underline N_\xi$
is of Kummer type.
\end{enumerate}
Then there exists a (Zariski) open neighborhood $U$ of\/ $|\xi|$ in
$Y$, such that the restriction $f_{|U}:(U,\underline
N_{|U})\to(X,\underline M)$ of $f$ is of Kummer type.
\end{lemma}
\begin{proof} By corollary \ref{cor_sharpy}(i) and theorem
\ref{th_good-charts}(ii), we may find a neighborhood $U'\to Y$ of
$\xi$ in $Y_\tau$, and a finite chart $(\omega_P,\omega_Q,\theta)$
for the restriction $f_{|U'}$, with $Q$ fine and saturated. Set
$$
S_P:=\omega_{P,\xi}^{-1}\underline N^\times_\xi \quad
S_Q:=\omega_{Q,f(\xi)}^{-1}\underline M^\times_{f(\xi)} \quad
P':=S_P^{-1}P \quad Q':=S_Q^{-1}Q.
$$
According to claim \ref{cl_better-chart}(iii) we may find
neighborhoods $U''\to U'$ of $\xi$ in $Y_\tau$, and $V\to X$
of $f(\xi)$ in $X_\tau$, such that the charts $\omega_{P|U''}$
and $\omega_{Q|V}$ extend to charts
$$
\omega_{P'}:P'_{U''}\to\underline N_{|U''} \qquad
\omega_{Q'}:Q'_{V}\to\underline M_{|V}.
$$
Clearly $\theta$ extends as well to a unique morphism $\theta':Q'\to
P'$, and after shrinking $U''$ we may assume that the restriction
$f_{|U''}:(U'',\underline N_{|U''})\to(X,\underline M)$ factors
through a morphism $f':(U'',\underline N_{|U''})\to(V,\underline
M_{|V})$, in which case it is easily seen that the datum
$(\omega_{P'},\omega_{Q'},\theta')$ is a chart for $f'$. Notice as
well that $Q'$ is still fine and saturated (lemma
\ref{lem_exc-satura}(i)), and by claim \ref{cl_better-chart}(iv) the
maps $\omega_{P'}$ and $\omega_{Q'}$ induce isomorphisms :
$$
P^{\prime\sharp}\isom\underline N{}^\sharp_\xi
\qquad
Q^{\prime\sharp}\isom\underline M{}_{f(\xi)}^\sharp.
$$
Our assumption (b) then implies that the map
$Q^{\prime\sharp}\to P^{\prime\sharp}$ deduced from $\theta'$
is of Kummer type, and then the morphism of fans
$\Spec\,\theta':\Spec\,P'\to\Spec\,Q'$ is of Kummer type as well
(remark \ref{rem_Kummer-type}(i)). However, let
$$
\bar\omega_{P'}:
(U'',\underline N{}^\sharp_{|U''})\to(\Spec\,P')^\sharp
\qquad 
\bar\omega_{Q'}:
(V,\underline M{}^\sharp_{|V})\to(\Spec\,Q')^\sharp
$$
be the morphisms of monoidal spaces deduced from $\omega_{P'}$ and
$\omega_{Q'}$; we obtain a morphism
$$
(f',\Spec\,\theta'):
(U'',\underline N_{|U''},(\Spec\,P')^\sharp,\bar\omega_{P'})\to
(V,\underline M_{|V},(\Spec\,Q')^\sharp,\bar\omega_{Q'})
$$
in the category $\cK$ of \eqref{subsec_category-K}, which -- in view
of remark \ref{rem_Kummer-type}(ii) -- shows that $f'$ is of Kummer
type. This already concludes the proof in case $\tau=\Zar$, and for
$\tau=\et$ it suffices to remark that the image of $U''$ in $Y$ is a
Zariski open neighborhood $U$ of $\xi$, such that the restriction
$f_{|U}$ is of Kummer type.
\end{proof}

We shall use the following criterion :

\begin{proposition}\label{prop_crit-Kummer}
Let $k$ be a field, $f:P\to Q$ an injective morphism of fine
monoids, with $P$ sharp. Suppose that the scheme $\Spec\,k\La
Q/\fm_P Q\Ra$ admits an irreducible component of Krull dimension
$0$. We have :
\begin{enumerate}
\item
$f$ is of Kummer type, and $k\La Q/\fm_P Q\Ra$ is a finite
$k$-algebra.
\item
If moreover, $Q^\times$ is a torsion-free abelian group, then $Q$ is
sharp, and $k\La Q/\fm_P Q\Ra$ is a local $k$-algebra with $k$ as
residue field.
\end{enumerate}
\end{proposition}
\begin{proof} Notice that the assumption means especially that
$Q/\fm_PQ\neq\{1\}$, so $f$ is a local morphism. Set
$I:=\rad(\fm_P Q)$ (notation of definition \ref{def_radical-mon}(ii)),
and let $\fp_1,\dots,\fp_n\subset Q$ be the minimal prime ideals
containing $I$; then $I=\fp_1\cap\cdots\cap\fp_n$, by lemma
\ref{lem_radical}.
Clearly the natural closed immersion
$\Spec\,k\La Q/I\Ra\to\Spec\,k\La Q/\fm_PQ\Ra$ is a homeomorphism;
on the other hand, say that $\fq\subset k\La Q\Ra$ is a prime ideal
containing $I$; then $\fq\cap Q$ is a prime ideal of $Q$ containing
$I$, hence $\fp_i\subset \fq$ for some $i\leq n$, {\em i.e.}
$\Spec\,k\La Q/I\Ra$ is the union of its closed subsets
$\Spec\,k\La Q/\fp_iQ\Ra$, for $i=1,\dots,n$. The Krull dimension
of each irreducible component of
$\Spec\,k\La Q/\fp_i\Ra=\Spec\,k[Q\!\setminus\!\fp_i]$ equals
$$
\rk_\Z(Q\!\setminus\!\fp)^\times+\dim Q\!\setminus\!\fp=
\rk_\Z Q^\times+d-\hgt(\fp_i)
\qquad
\text{where $d:=\dim Q$}
$$
(claim \ref{cl_NullSt}(ii) and corollary \ref{cor_consequent}(i,ii)).
Our assumption then implies that $Q^\times$ is a finite group,
and $\hgt(\fp_i)=d$, {\em i.e.} $\fp_i=\fm_Q$, for at least an
index $i\leq n$, and therefore $I=\fm_Q$.
Furthermore, since $\fm_Q$ is finitely generated, we have
$\fm^n_Q\subset\fm_PQ$ for a sufficiently large integer $n>0$, and
then it follows easily that $k\La Q/\fm_P Q\Ra$ is a finite
$k$-algebra. If $Q^\times$ is torsion-free, then we also deduce
that $Q$ is sharp, and moreover the maximal ideal of
$k\La Q/\fm_P Q\Ra$ generated by $\fm_Q$ is nilpotent, hence the
latter $k$-algebra is local.

Let now $\fq\subset Q$ be any prime ideal, and pick
$x\in\fm_Q\!\setminus\!\fq$; the foregoing implies that there exists
an integer $r>0$ such that $x^r=f(p_\fq)q$ for some $p_\fq\in\fm_P$
and $q\in Q$, hence $f(p_\fq)\notin\fq$, and $f(p_\fq)$ is not
invertible in $Q$, since $f$ is local. If now $\fq$ has height
$d-1$, it follows that $f(p_\fq)$ is a generator of
$(Q\!\setminus\!\fq)_\R$, which is an extremal ray of the polyhedral
cone $Q_\R$, and every such extremal ray is of this form
(proposition \ref{prop_Gordie}); furthermore, the latter cone is
strictly convex, since $Q^\times$ is finite. Hence the set
$S:=(f(p_\fq)~|~\hgt(\fq)=d-1)$ is a system of generators of the
polyhedral cone $Q_\R$ (see \eqref{subsec_extremal}). Let $Q'\subset
Q$ be the submonoid generated by $S$; then $S_\Q=S_\R\cap Q_\Q=Q_\Q$
(proposition \ref{prop_Gordon}(iii)), {\em i.e.} $f$ is of Kummer
type.
\end{proof}

\sset\subsubsection{}
Let $f:(Y_\Zar,\underline N)\to(X_\Zar,\underline M)$ be a morphism
of log schemes with Zariski log structures; it follows easily from
the isomorphism \eqref{eq_general-iso} and remark
\ref{rem_Kummer-type}(ii) that $f$ is of Kummer type if and only if
$\tilde u{}^*f:(X_\et,\tilde u{}^*\underline M)\to
(Y_\et,\tilde u{}^*\underline N)$ is a morphism of Kummer type
between schemes with \'etale log structures (notation of
\eqref{subsec_choose-a-top}). Suppose now that $\underline M$ is
an integral and saturated log structure on $X_\Zar$, and denote by
$\mathbf{s.Kum}(X_\Zar,\underline M)$ (resp.
$\mathbf{s.Kum}(X_\et,\tilde u{}^*\underline M)$) the full
subcategory of $\bsatlog/(X_\Zar,\underline M)$ (resp. of
$\bsatlog/(X_\et,\tilde u{}^*\underline M)$ whose objects are all
the morphisms of Kummer type. In view of the foregoing (and of
lemma \ref{lem_simple-charts-top}(i)), we see that $\tilde u{}^*$
restricts to a functor :
\set\begin{equation}\label{eq_pulbak-kummer}
\mathbf{s.Kum}(X_\Zar,\underline M)\to\mathbf{s.Kum}(X_\et,\tilde
u{}^*\underline M).
\end{equation}

\begin{lemma}\label{lem_Kummer-descends}
The functor \eqref{eq_pulbak-kummer} is an equivalence.
\end{lemma}
\begin{proof} By virtue of proposition
\ref{prop_reduce-to-etale}(ii) we know already that
\eqref{eq_pulbak-kummer} is fully faithful, hence we only
need to show its essential surjectivity. Thus, let
$f:(Y_\et,\underline N)\to(X_\et,\tilde u{}^*\underline M)$
be a morphism of Kummer type. By remark \ref{rem_Kummer-type}(ii),
we know that $\log f$ induces an isomorphism
$$
\tilde u{}^*f^*\underline M{}^\sharp_\Q\isom
f^*\tilde u{}^*\underline M{}^\sharp_\Q\isom
\underline N{}^\sharp_\Q
$$
(notation of \eqref{subsec_from-con-to-mon}). Since $\underline N$
is integral and saturated, the natural map
$\underline N^\sharp\to\underline N{}^\sharp_\Q$ is a monomorphism,
so that the counit of adjunction
$\tilde u{}^*\tilde u_*\underline N^\sharp\to\underline N^\sharp$
is an isomorphism (lemma \ref{lem_Hilbert90}(ii)), and then the
same holds for the counit of adjunction
$\tilde u{}^*\tilde u_*(Y,\underline N)\to(Y,\underline N)$
(proposition \ref{prop_reduce-to-etale}(iii)).
\end{proof}

\begin{proposition}\label{prop_characterize-log-etale-cov}
Let $f:(Y,\underline N)\to(X,\underline M)$ be a morphism of fs log
schemes. The following conditions are equivalent:
\begin{enumerate}
\alphaenu
\item
Every geometric point of $X$ admits an \'etale neighborhood $U\to X$
such that $Y_U:=U\times_XY$ decomposes as a disjoint union
$Y_U=\bigcup_{i=1}^nY_i$ of open and closed subschemes (for some
$n\in\N$), and we have :
\begin{enumerate}
\romanenuii
\item
Each restriction $Y_i\times_Y(Y,\underline N)\to
U\times_X(X,\underline M)$ of $f\times_X\one_U$ admits a fine,
saturated Kummer chart $(\omega_{P_i},\omega_{Q_i},\theta_i)$ such
that $P_i$ and $Q_i$ are sharp, and the order of\/
$\Coker\,\theta^\gp_i$ is invertible in $\cO_{\!U}$.
\item
The induced morphism of\/ $U$-schemes $Y_i\to
U\times_{\Spec\,\Z[P_i]}\Spec\,\Z[Q_i]$ is an isomorphism, for every
$i=1,\dots,n$.
\end{enumerate}
\item
$f$ is \'etale, and the morphism of schemes underlying $f$ is
finite.
\end{enumerate}
\end{proposition}
\begin{proof}(a) $\Rightarrow$ (b): Indeed, it is easily seen that
the morphism $\Spec\,\Z[\theta_i]$ is finite, hence the same holds
for the restriction $Y_i\to U$ of $f$, in view of (a.ii), and then
the same holds for $f\times_X\one_U$, so finally $f$ is finite on
the underlying schemes (\cite[Ch.IV, Prop.2.7.1]{EGAIV-2}), and it
is \'etale by the criterion of theorem \ref{th_charact-smoothness}.

(b) $\Rightarrow$ (a): Arguing as in the proof of theorem
\ref{th_charact-smoothness}, we may reduce to the case where
$\tau=\et$. Suppose first that both $X$ and $Y$ are strictly local.
In this case, $\underline M$ admits a fine and saturated chart
$\omega_P:P_X\to\underline M$, sharp at the closed point (corollary
\ref{cor_sharpy}(i)). Moreover, $f$ admits a chart
$(\omega_P,\omega_Q:Q_Y\to\underline N,\theta:P\to Q)$, for some
fine monoid $Q$ such that $Q^\times$ is torsion-free; also $\theta$
is injective, the induced morphism of $X$-schemes
$$
g:Y\to X':=X\times_{\Spec\,\Z[P]}\Spec\,\Z[Q]
$$
is \'etale, and the order of $\Coker\,\theta^\gp$ is invertible in
$\cO_{\!X}$ (corollary \ref{cor_charact-smoothness}). Since $f$ is
finite, $g$ is also closed (\cite[Ch.II, Prop.6.1.10]{EGAII}), hence
its image $Z$ is an open and closed local subscheme, finite over
$X$. Let $k$ be the residue field of the closed point of $X$; it
follows that $X'\times_X\Spec\,k\simeq\Spec\,k\La Q/\fm_P Q\Ra$
admits an irreducible component of Krull dimension zero (namely, the
intersection of $Z$ with the fibre of $X'$ over the closed point of
$X$), in which case the criterion of proposition
\ref{prop_crit-Kummer}(ii) ensures that $\theta$ is of Kummer type
and $Q$ is sharp, hence $X'$ is finite over $X$, and moreover
$X'\times_X\Spec\,k$ is a local scheme with $k$ as residue field.
Since $X$ is strictly henselian, it follows that $X'$ itself is
strictly local, and therefore $g$ is an isomorphism, so the
proposition is proved in this case.

Let now $X$ be a general scheme, and $\xi$ a geometric point of $X$;
denote by $X(\xi)$ the strict henselization of $X$ at the point
$\xi$, and set $Y(\xi):=X(\xi)\times_XY$. Since $f$ is finite,
$Y(\xi)$ decomposes as the disjoint union of finitely many open and
closed strictly local subschemes $Y_1(\xi),\dots,Y_n(\xi)$. Then we
may find an \'etale neighborhood $U\to X$ of $\xi$, and open and
closed subschemes $Y_1,\dots,Y_n$ of $Y\times_XU$, with isomorphisms
of $X(\xi)$-schemes $Y_i(\xi)\isom Y_i\times_UX(\xi)$, for every
$i=1,\dots,n$ (\cite[Ch.IV, Cor.8.3.12]{EGAIV-3}). We may then
replace $X$ by $U$, and we reduce to proving the proposition for
each of the restrictions $Y_i\to U$ of $f\times_X\one_U$; hence we
may assume that $Y(\xi)$ is strictly local. By the foregoing case,
we may find a chart
\set\begin{equation}\label{eq_given-chart}
P_{X(\xi)}\to\underline M(\xi) \qquad Q_{Y(\xi)}\to\underline N(\xi)
\qquad \theta:P\to Q
\end{equation}
of $f\times_X\one_{X(\xi)}$, with $\theta$ of Kummer type, such that
$P$ and $Q$ are sharp, the order $d$ of $\Coker\,\theta^\gp$ is
invertible in $\cO_{\!X(\xi)}$, and the induced morphism of
$X(\xi)$-schemes $Y(\xi)\to X(\xi)\times_{\Spec\,\Z[P]}\Spec\,\Z[Q]$
is an isomorphism.
By corollary \ref{cor_descend-chart-from-infty} we may
find an \'etale neighborhood $U\to X$ of $\xi$ such that
\eqref{eq_given-chart} extends to a chart for $f\times_X\one_U$.
After shrinking $U$, we may assume that $d$ is invertible in
$\cO_{\!U}$. Lastly, after further shrinking of $U$, we may
ensure that the induced morphism of $U$-schemes
$U\times_XY\to U\times_{\Spec\,\Z[P]}\Spec\,\Z[Q]$ is an
isomorphism (\cite[Ch.IV, Cor.8.8.2.4]{EGAIV-3}).
\end{proof}

\begin{definition} Let $f:(Y,\underline N)\to(X,\underline M)$ be a
morphism of fs log schemes. We say that $f$ is an {\em \'etale
covering\/} of $(X,\underline M)$, if $f$ fulfills the equivalent
conditions (a) and (b) of proposition
\ref{prop_characterize-log-etale-cov}. We denote by
$$
\bCov(X,\underline M)
$$
the full subcategory of the category of $(X,\underline M)$-schemes,
whose objects are the \'etale coverings of $(X,\underline M)$.
\end{definition}

\begin{remark}
(i)\ \
Notice that all morphisms in $\bCov(X,\underline M)$ are \'etale
coverings, in light of corollary \ref{cor_sorite-smooth}(ii).

(ii)\ \
Moreover, $\bCov(X,\underline M)$ is a Galois category
(see \cite[Exp.V, D\'ef.5.1]{SGA1}), and if $\xi$ is any
geometric point of $(X,\underline M)_\tr$, we obtain a
fibre functor for this category, by the rule :
$f\mapsto f^{-1}(\xi)$, for every \'etale covering $f$
of $(X,\underline M)$. (Details left to the reader.)
We shall denote by
$$
\pi_1((X,\underline M)_\et,\xi)
$$
the corresponding fundamental group.
\end{remark}

\begin{example}\label{ex_Kummer-is-tamely-ram}
Let $(f,\log f):(Y,\underline N)\to(X,\underline M)$ be an \'etale
covering of a regular log scheme $(X,\underline M)$, and suppose
that $X$ is strictly local of dimension $1$ and $Y$ is connected
(hence strictly local as well). Let $x\in X$ (resp. $y\in Y$) be the
closed point; it follows that $\cO_{\!X,x}$ is a strictly henselian
discrete valuation ring (corollary \ref{cor_normal-and-CM}). The
same holds for $\cO_{\!Y,y}$ in view of theorem
\ref{th_smooth-preserve-reg} and \cite[Ch.IV, Prop.18.5.10]{EGA4}.
In case $\dim\underline M_x=0$, the log structures $\underline M$
and $\underline N$ are trivial, so $f:Y\to X$ is an \'etale morphism
of schemes. Otherwise we have $\dim\underline M_x=1=\dim\underline N_y$
(lemma \ref{lem_Kummer-fans}(i)), and then
$\underline M{}_x^\sharp\simeq\N\simeq\underline N{}_y^\sharp$
(theorem \ref{th_structure-of-satu}(ii)); also, the choice of a chart
for $f$ as in proposition \ref{prop_characterize-log-etale-cov},
induces an isomorphism
$$
\cO_{\!X,x}\otimes_{\underline M{}_x^\sharp}\underline N{}^\sharp_y
\isom\cO_{Y,y}
$$
where the map $\underline M{}_x^\sharp\to\underline N{}^\sharp_y$
is the $N$-Frobenius map of $\N$, where $N>0$ is an integer
invertible in $\cO_{\!X,x}$. Moreover, notice that the image
of the maximal ideal of $\underline M_x$ generates the maximal
ideal of $\cO_{\!X,x}$ (and likewise for the image of
$\underline N{}_y$ in $\cO_{Y,y}$), so the structure map of
$\underline M$ induces an isomorphism
$\underline M{}_x^\sharp\isom\Gamma_{\!+}$, onto the submonoid of
the value group $(\Gamma,\leq)$ of $\cO_{\!X,x}$ consisting of
all elements $\leq 1$ (and likewise for $\underline N{}_y^\sharp$).
In other words, $\cO_{Y,y}$ is obtained from $\cO_{X,x}$ by
adding the $N$-th root of a uniformizer.
It then follows that the ring homomorphism $\cO_{\!X,x}\to\cO_{Y,y}$
is an algebraic {\em tamely ramified extension\/} of discrete
valuation rings (see {\em e.g.} \cite[Cor.6.2.14]{Ga-Ra}).
\end{example}

\begin{definition}\label{def_tamely-ram}
Let $X$ be a normal scheme, and $Z\to X$ a closed immersion such
that :
\begin{enumerate}
\alphaenu
\item
$Z$ is a union of irreducible closed subsets of codimension one in $X$.
\item
For every maximal point $z\in Z$, the stalk $\cO_{\!X,z}$ is a
discrete valuation ring.
\end{enumerate}
Set $U:=X\!\setminus\!Z$, let $f:U'\to U$ be an \'etale covering,
and $g:X'\to X$ the normalization of $X$ in $U'$. Notice that
$g^{-1}Z$ is a union of irreducible closed subsets of codimension
one in $X'$ (by the going down theorem \cite[Th.9.4(ii)]{Mat}).
\begin{enumerate}
\item
We say that $f$ is {\em tamely ramified along $Z$} if, for every
geometric point $\xi\in X'$ localized at a maximal point of
$g^{-1}Z$, the induced finite extension
$$
\cO^\sh_{\!X,g(\xi)}\to\cO^\sh_{\!X',\xi}
$$
of strictly henselian discrete valuation rings
(\cite[Ch.IV, Cor.18.8.13]{EGA4}), is tamely ramified. Recall that
the latter condition means the following  :  let
$\Gamma_{\!\!g(\xi)}\to\Gamma_{\!\xi}$ be the map of value
groups induced by the above extension; then the
{\em ramification index of $f$ at $\xi$}
$$
e_\xi(f):=(\Gamma_{\!\xi}:\Gamma_{\!\!f(\xi)})
$$
is invertible in the residue field $\kappa(\xi)$, and $g$ induces
an isomorphism $\kappa(f(\xi))\isom\kappa(\xi)$.
\item
Suppose that $f$ is tamely ramified along $Z$, and that the set
of maximal points of $Z$ is finite. Then the {\em ramification
index of $f$ along $Z$\/} is the least common multiple $e_Z(f)$
of the ramification indices $e_\xi(f)$, where $\xi$ ranges over
the set of all geometric points of $X'$ supported at a maximal
point of $g^{-1}Z$.
\item
We denote by $\bTame(X,U)$ the full subcategory of $\bCov(U)$
whose objects are the \'etale coverings of $U$ that are
tamely ramified along $Z$.
\end{enumerate}
\end{definition}

\begin{lemma}\label{lem_tameness-pullsback}
Let $\phi:X'\to X$ be a dominant morphism of normal schemes,
$Z\subset X$ a closed subset, and set $Z':=\phi^{-1}Z$. Suppose
that :
\begin{enumerate}
\alphaenu
\item
The closed immersions $Z\to X$ and $Z'\to X'$ satisfy conditions
{\em (a)} and {\em (b)} of definition {\em\ref{def_tamely-ram}}.
\item
$\phi$ restricts to a map $\Max\,Z'\to\Max\,Z$ on the subsets
of maximal points of $Z$ and $Z'$.
\end{enumerate}
We have :
\begin{enumerate}
\item
$\phi^*:\bCov(U)\to\bCov(U')$ restricts to a functor
$$
\bTame(X,U)\to\bTame(X',U').
$$
\item
Suppose that $\Max\,Z$ and $\Max\,Z'$ are finite sets. Then, for
any object $f:V\to U$ of\/ $\bTame(X,U)$, the index $e_{Z'}(\phi^*f)$
divides $e_Z(f)$.
\item
Suppose that $X$ is regular, and let $j:U\to X$ be the open immersion.
Then the essential image of the pull-back functor
$j^*:\bCov(X)\to\bTame(X,U)$ consists of the objects $f:V\to U$ such
that $e_Z(f)=1$.
\end{enumerate}
\end{lemma}
\begin{proof} Let $V\to U$ be an object of $\bTame(X,U)$, and
denote by $g:W\to X$ (resp. $g':W'\to X'$) the normalization of $X$
in $V$ (resp. of $X'$ in $V\times_UU'$). Let $w'$ be a geometric
point of $W'$ localized at a maximal point of $g^{\prime-1}Z'$,
and denote by $w$ (resp. $x'$, resp. $x$) the image of $w'$ in
$W$ (resp. in $X'$, resp. in $X$); in order to check (i) and (ii),
we have to show that the ramification index of $g'$ at the
geometric point $w'$ is invertible in $\kappa(w')$, and the
residue field extension $\kappa(x')\subset\kappa(w')$ is trivial.
To this aim, in view of \cite[Ch.IV, Prop.17.5.7]{EGA4}, we may
replace $X'$ by $X'(x')$, $X$ by $X(x)$, $W$ by $W(w)$, and
assume from start that $X$, $X'$ and $W$ are strictly local,
say $X=\Spec\,A$, $X'=\Spec\,B$ and $C=\Spec\,W$ for some
strictly henselian local rings $A,B,C$. In this case,
$D:=B\otimes_AC$ is a direct product of finite and local
$B$-algebras, and $x'$ is localized at a closed point of
$X'$, so $w'$ is localized at a closed point of $W'$,
therefore $\cO_{W',w'}$ is a direct factor of $D$, and
$\kappa(w')$ is a quotient of
$\kappa(x')\otimes_{\kappa(x)}\kappa(w)$; but we have
$\kappa(w)=\kappa(x)$ by assumption, so this already
yields $\kappa(w')=\kappa(x')$. Moreover,
by construction, $\cO_{W',w'}\otimes_A\Frac\,A$
is a finite separable extension of $\Frac\,B$, and the
diagram
\set\begin{equation}\label{eq_finite-ab-vf-exts}
{\diagram
\Frac\,A \ar[r] \ar[d] & \Frac\,C \ar[d] \\
\Frac\,B \ar[r] & \Frac\,\cO_{W',w'}
\enddiagram}
\end{equation}
is cocartesian (in the category of fields). Since $V$ is tamely
ramified along $Z$, the top horizontal arrow of
\eqref{eq_finite-ab-vf-exts} is a finite (abelian) extension whose
degree $e_w$ is invertible in $\kappa(x)$
(\cite[Cor.6.2.14]{Ga-Ra}); it follows that the bottom horizontal
arrows is a Galois extension as well, and its degree divides $e_w$,
as required.

(iii): Let $f:V\to U$ be an object of $\bTame(X,U)$.
By claim \ref{cl_trivial-Zorni}, there exists a largest open subset
$U_{\max}\subset X$ containing $U$, and such that $f$ is the
restriction of an \'etale covering $f_{\max}:V'\to U_{\max}$;
now, if $e_Z(f)=1$, every point of codimension one of $X\!\setminus\!U$
lies in $U_{\max}$, by claim \ref{cl_include-pts}. Hence $X\!\setminus\!U$
has codimension $\geq 2$ in $X$, in which case lemma \ref{lem_replace}(iv)
implies that $X=U_{\max}$, as required.
\end{proof}

\begin{remark} The assumptions (a) and (b) of lemma
\ref{lem_tameness-pullsback} are fulfilled, notably, when
$f$ is finite and dominant, or when $f$ is flat
(\cite[Th.9.4(ii), Th.9.5]{Mat}).
\end{remark}

\sset\subsubsection{}\label{subsec_many-condits}
Let now $\underline X:=(X_i~|~i\in I)$ be a cofiltered system
of normal schemes, such that, for every morphism $\phi:i\to j$
in the indexing category $I$, the corresponding transition
morphism $f_\phi:X_i\to X_j$ is dominant and affine. Suppose
also, that for every $i\in I$, there exists a closed immersion
$Z_i\to X_i$, fulfilling conditions (a) and (b) of definition
\ref{def_tamely-ram}, such that $\Max\,Z_i$ is a finite set,
and for every morphism $\phi:i\to j$ of $I$, we have :
\begin{itemize}
\item
$Z_i=f_\phi^{-1}Z_j$.
\item
The corresponding morphism $f_\phi$ restricts to a map
$\Max\,Z_i\to\Max\,Z_j$.
\end{itemize}
Let also $X$ be the limit of $\underline X$, and $Z$ the
limit of the system $(Z_i~|~i\in I)$, and suppose furthermore
that the closed immersion $Z\to X$ fulfills as well conditions
(a) and (b) of definition \ref{def_tamely-ram}.

\begin{proposition}\label{prop_tame-infty}
In the situation of \eqref{subsec_many-condits}, we have :
\begin{enumerate}
\item
The induced morphism $\pi_i:X\to X_i$ restricts to a map
$$
\Max\,Z\to\Max\,Z_i
\qquad
\text{for every $i\in I$.}
$$
\item
The morphisms $\pi_i$ induce an equivalence of categories :
$$
\Pscolim{I}\bTame(X_i,X_i\!\setminus\!Z_i)\to\bTame(X,X\!\setminus\!Z).
$$
\end{enumerate}
\end{proposition}
\begin{proof}(i): More precisely, we shall prove that there is a
natural homeomorphism :
$$
\Max\,Z\isom\lim_{i\in I}\Max\,Z_i.
$$
(Notice the each $\Max\,Z_i$ is a discrete finite set, hence this
will show that $\Max\,Z$ is a profinite topological space.)
Indeed, suppose that $z:=(z_i~|~i\in I)$ is a maximal point of $Z$.
For every $i\in I$, let $T_i\subset\Max\,Z_i$ be the subset
of maximal generizations of $z_i$ in $Z_i$. It is easily seen
that $f_\phi T_i\subset T_j$ for every morphism $\phi:i\to j$
in $I$. Clearly $T_i$ is a finite non-empty set for every
$i\in I$, hence the limit $T$ of the cofiltered system
$(T_i~|~i\in I)$ is non-empty. However, any point of $T$
is a generization of $z$ in $Z$, hence it must coincide with
$z$. The assertion follows easily.

(ii): To begin with, lemma \ref{lem_tameness-pullsback}(i)
shows that, for every morphism $\phi:i\to j$ in $I$, the
transition morphism $f_\phi$ induces a pull-back functor
$\bTame(X_j,X_j\!\setminus\!Z_j)\to\bTame(X_i,X_i\!\setminus\!Z_i)$,
so the $2$-colimit in (ii) is well-defined, and combining (i)
with lemma \ref{lem_tameness-pullsback}(i) we obtain indeed
a well-defined functor from this $2$-colimit to
$\bTame(X,X\!\setminus\!Z)$.

The full faithfulness of the functor of (ii) follows
from lemma \ref{lem_go-to-pseudo-lim}. Next, let
$g:V\to X\!\setminus\!Z$ be an object of $\bTame(X,X\!\setminus\!Z)$;
invoking again lemma \ref{lem_go-to-pseudo-lim}, we may
descend $g$ to an \'etale covering $g_j:V_j\to X_j\!\setminus\!Z_j$,
for some $j\in I$; after replacing $I$ by $I/j$, we may assume that
$j$ is the final object of $I$ and we may define $g_i:=f^*_\phi(g_j)$
for every $\phi:i\to j$ in $I$. To conclude the proof, it suffices to
show that there exists $i\in I$ such that $g_i$ is tamely ramified
along $Z_i$.

Now, let $\bar g_i:\bar V_i\to X_i$ be the normalization
of $X_i$ in $V_i\times_{X_j}X_i$, for every $i\in I$, and
$\bar g:\bar V\to X$ the normalization of $X$ in $V$.
Given a geometric point $\bar v$ localized at a maximal point
of $\bar g{}^{-1}Z$, let $\bar v_i$ (resp. $\bar z_i$) be the
image of $\bar v$ in $\bar V_i$ (resp. in $Z_i$), for every
$i\in I$. Let also $\bar z$ be the image of $\bar v$ in $Z$.
Then
$$
\cO^\sh_{\!X,\bar z}=\colim_{i\in I}\cO^\sh_{\!X_i,\bar z_i}
\qquad
\cO^\sh_{\bar V,\bar v}=\colim_{i\in I}\cO^\sh_{\bar\V_i,\bar v_i}.
$$
(\cite[Ch.IV, Prop.18.8.18(ii)]{EGA4}), and it follows easily
that there exists $i\in I$ such that the finite extension
$\cO^\sh_{\!X_i,\bar z_i}\to\cO^\sh_{\bar\V_i,\bar v_i}$ is
already tamely ramified. Since only finitely many points of
$\bar V_i$ lie over the support of $\bar z_i$, and since $I$
is cofiltered, it follows that there exists $i\in I$ such that
the induced morphism
$\bar V_i\times_{X_i}X_i(\bar z_i)\to X_i(\bar z_i)$
is already tamely ramified. However, notice that $\pi^{-1}_i(z_i)$
is open in $\Max\,Z$, for every $i\in I$, and every $z_i\in\Max\,Z_i$.
Therefore, we may find a finite subset $I_0\subset I$, and for
every $i\in I_0$ a subset $T_i\subset\Max\,Z_i$, such that :
\begin{itemize}
\item
$\Max\,Z=\bigcup_{i\in I_0}\pi_i^{-1}(T_i)$.
\item
For every geometric point $\bar z_i$ localized in $T_i$,
the morphism
$\bar V_i\times_{X_i}X_i(\bar z_i)\to X_i(\bar z_i)$
is tamely ramified.
\end{itemize}
Since $I$ is cofiltered, we may find $k\in I$ with morphisms
$\phi_i:k\to i$ for every $i\in I_0$; after replacing $T_i$
by $f_{\phi_i}^{-1}(T_i)$ for every $i\in I_0$, we may then
assume that $I_0=\{k\}$, so that $\Max\,Z=\pi_k^{-1}(T_k)$.
It follows that $\Max\,Z_i=f_\phi^{-1}T_k$ for some $i\in I$
and some $\phi:i\to k$. Then it is clear that $g_i$ is tamely
ramified along $Z_i$, as required.
\end{proof}

\sset\subsubsection{}\label{subsec_Kummer-exact-seq}
Let $X$ be a scheme, $U\subset X$ a connected open subset, $\xi$ any
geometric point of $U$, and $N>0$ an integer which is invertible in
$\cO_{\!U}$; then the $N$-Frobenius map $\boldsymbol{N}$ of
$\cO^\times_{\!U}$ gives a {\em Kummer exact sequence\/} of
abelian sheaves on $X_\et$ :
$$
0\to\bmu_{N,U}\to\cO^\times_{\!U}\xrightarrow{\ \boldsymbol{N}\ }
\cO^\times_{\!U}\to 0
$$
(where $\bmu_{N,U}$ is the $N$-torsion subsheaf of
$\cO^\times_{\!U_\et}$ : see \cite[Exp.IX, \S3.2]{SGA4-3}). Suppose
now that $\bmu_{N,U}$ is a constant sheaf on $U_\et$; this means
especially that $\bmu_N:=(\bmu_{N,U})_\xi$ is a cyclic group of
order $N$. Indeed, $\bmu_N$ certainly contains such a subgroup
(since $N$ is invertible in $\cO_{\!U}$), so denote by $\zeta$
one of its generators; then every $u\in\bmu_N$ satisfies the
identity $0=u^N-1=\prod_{i=1}^N(u-\zeta^i)$, and each factor
$u_\zeta^i$ of this decomposition vanishes on a closed subset
$U_i$ of $U(\xi)$; clearly $U_i\cap U_j=\emptyset$ for $i\neq j$,
so we get a decomposition of $U(\xi)$ as a disjoint union of open
and closed subsets $U(\xi)=U_1\cup\cdots\cup U_N$ such that
$u=\zeta^i$ on $U_i$, for every $i\leq N$; since $U(\xi)$ is
connected, it follows that $U(\xi)=U_i$ for some $i$, {\em i.e.}
$\zeta$ generates $\bmu_N$. There follows a natural map :
$$
\partial_N:\Gamma(U_\et,\cO^\times_{\!U})\to H^1(U_\et,\bmu_{N,U})\isom
\Hom_\mathrm{cont}(\pi_1(U_\et,\xi),\bmu_N)
$$
(lemma \ref{lem_from-coh-to-fundgrp}). Recall the geometric
interpretation of $\partial_N$ : a given $u\in\cO^\times_{\!U}(U)$
is viewed as a morphism of schemes $u:U\to\G_m$, where
$\G_m:=\Spec\,\Z[\Z]$ is the standard multiplicative group
scheme; we let $U'$ be the fibre product in the cartesian diagram
of schemes :
\set\begin{equation}\label{eq_G_m-torsor}
{\diagram
U' \ar[r] \ar[d]_{\phi_u} &
\G_m \ar[d]^{\Spec\,\Z[\boldsymbol{N}_\Z]} \\
U \ar[r]^-u & \G_m.
\enddiagram}
\end{equation}
Then $\phi_u$ is a torsor under the $U_\et$-group $\bmu_{N,U}$, and
to such torsor, lemma \ref{lem_from-coh-to-fundgrp} attaches a
well defined continuous group homomorphism as required.

If $M>0$ is any integer dividing $N$, a simple inspection yields a
commutative diagram :
\set\begin{equation}\label{eq_assembly-lines}
{\diagram \Gamma(U_\et,\cO^\times_{\!X}) \ar[rr]^-{\partial_N}
\ar[rrd]_{\partial_M} & &
\Hom_\mathrm{cont}(\pi_1(U_\et,\xi),\bmu_N) \ar[d]^{\pi_{N,M}} \\
& & \Hom_\mathrm{cont}(\pi_1(U_\et,\xi),\bmu_M) \enddiagram}
\end{equation}
where $\pi_{N,M}$ is induced by the map $\bmu_N\to\bmu_M$ given
by the rule : $x\mapsto x^{N/M}$ for every $x\in\bmu_N$.

\sset\subsubsection{}\label{subsec_first-pairing}
Let now $X$ be a strictly local scheme, $x$ the closed point
of $X$, and denote by $p$ the characteristic exponent of the
residue field $\kappa(x)$ (so $p$ is either $1$ or a positive
prime integer). Let also $\beta:\underline M\to\cO_{\!X}$ be a
log structure on $X_\et$, take $U:=(X,\underline M)_\tr$,
suppose that $U\neq\emptyset$, and fix a geometric point $\xi$
of $U$; for every integer $N>0$ with $(N,p)=1$, we get a morphism
of monoids :
\set\begin{equation}\label{eq_getting-there}
\underline M{}_x\xrightarrow{\ \beta_x\ }
\Gamma(U_\et,\cO^\times_{\!X})\xrightarrow{\ \partial_N\ }
\Hom_\mathrm{cont}(\pi_1(U_\et,\xi),\bmu_N)
\end{equation}
whose kernel contains $\underline M{}_x^N$, the image of the
$N$-Frobenius endomorphism of $\underline M{}_x$. Notice that the
$N$-Frobenius map of $\cO^\times_{\!X,x}$ is surjective : indeed,
since the residue field $\kappa(x)$ is separably closed, and
$(N,p)=1$, all polynomials in $\kappa(x)[T]$ of the form $T^N-u$
(for $u\neq 0$) split as a product of distinct monic polynomials
of degree $1$; since $\cO_{\!X,x}$ is henselian, the same holds
for all polynomials in $\cO_{\!X,x}[T]$ of the form $T^N-u$, with
$u\in\cO^\times_{\!X,x}$. It follows that \eqref{eq_getting-there}
factors through a natural map:
$$
\underline M{}^{\sharp\gp}_x\to
\Hom_\mathrm{cont}(\pi_1(U_\et,\xi),\bmu_N)
$$
which is the same as a group homomorphism :
\set\begin{equation}\label{eq_first-pairing}
\underline M{}^{\sharp\gp}_x\times\pi_1(U_\et,\xi)\to\bmu_N.
\end{equation}
In view of the commutative diagram \eqref{eq_assembly-lines},
it is easily seen that the pairings \eqref{eq_first-pairing},
for $N$ ranging over all the positive integers with $(N,p)=1$,
assemble into a single pairing :
$$
\underline M{}_x^\gp\times\pi_1((X,\underline M)_{\tr,\et},\xi)\to
\prod_{\ell\neq p}\Z_\ell(1)
$$
(where $\ell$ ranges over the prime numbers different from $p$, and
$\Z_\ell(1):=\lim_{n\in\N}\,\bmu_{\ell^n}$). The latter is the
same as a group homomorphism :
\set\begin{equation}\label{eq_tame-quotient}
\pi_1((X,\underline M)_{\tr,\et},\xi)\to
\underline M{}_x^{\gp\vee}\otimes_\Z\prod_{\ell\neq p}\Z_\ell(1).
\end{equation}

\sset\subsubsection{}\label{subsec_funct-standard-shit}
Let $f:(X,\underline M)\to(Y,\underline N)$ be a morphism of
log schemes, such that both $X$ and $Y$ are strictly local,
and $f$ maps the closed point $x$ of $X$ to the closed point
$y$ of $Y$. Then $f$ restricts to a morphism of schemes
$f_\tr:(X,\underline M)_\tr\to(Y,\underline N)_\tr$ (remark
\ref{rem_trivial-is-strict}(i)). Fix again a geometric point
$\xi$ of $(X,\underline M)_\tr$; we get a diagram of group
homomorphisms :
$$
\xymatrix{
\pi_1((X,\underline M)_{\tr,\et},\xi)
\ar[rr]^-{\pi_1(f_\tr,\xi)} \ar[d]
& & \pi_1((Y,\underline N)_{\tr,\et},f(\xi)) \ar[d] \\
\underline N{}^{\gp\vee}_y\otimes_\Z\prod_{\ell\neq p}\Z_\ell(1)
\ar[rr] & &
\underline M{}^{\gp\vee}_x\otimes_\Z\prod_{\ell\neq p}\Z_\ell(1)
}$$
whose vertical arrows (resp. bottom horizontal arrow) are the
maps \eqref{eq_tame-quotient} (resp. is induced by
$(\log f_x)^{\gp\vee}$), and by inspecting the constructions,
it is easily seen that this diagram commutes.

\sset\subsubsection{}\label{subsec_Second-pairing}
Suppose now that $(X,\underline M)$ as in
\eqref{subsec_first-pairing} is a fs log scheme, so that
$P:=\underline M{}^\sharp_x$ is fine and saturated; in this
case, we wish to give a second construction of the pairing
\eqref{eq_first-pairing}.

$\bullet$\ \
Namely, for every integer $N>0$
set $S_P:=\Spec\,\Z[P]$ and consider the finite morphism
of schemes $g_N:S_P\to S_P$ induced by the $N$-Frobenius
endomorphism of $P$. Notice that $S_P$ contains the dense
open subset $U_P:=\Spec\,\Z[1/N,P^\gp]$, and it is easily
seen that the restriction $g_{N|U_P}:U_P\to U_P$ of $g_N$
is an \'etale morphism.
Let $\tau\in U_P$ be any geometric point; then $\tau$ corresponds to
a ring homomorphism $\Z[1/N,P^\gp]\to\kappa:=\kappa(\tau)$, which is
the same as a character $\chi_\tau:P^\gp\to\kappa^\times$. Likewise,
any $\tau'\in g_N^{-1}(\tau)$ is determined by a character
$\chi_{\tau'}:P^\gp\to\kappa^\times$ extending $\chi_\tau$, {\em
i.e.} such that $\chi_{\tau'}(x^N)=\chi_\tau(x)$ for every $x\in
P^\gp$. Notice that every character $\chi_\tau$ as above admits at
least one such extension $\chi_{\tau'}$, since $\kappa$ is separably
closed, and $N$ is invertible in $\kappa$. Let $C_N$ be the cokernel
of the $N$-Frobenius endomorphism of $P^\gp$; there follows a short
exact sequence of finite abelian groups :
$$
0\to\Hom_\Z(C_N,\kappa^\times)\to\Hom_\Z(P^\gp,\kappa^\times)
\xrightarrow{ N }\Hom_\Z(P^\gp,\kappa^\times)\to 0.
$$
Especially, we see that the fibre $g_N^{-1}(\tau)$ is naturally a
torsor under the group
$\Hom_\Z(C_N,\kappa^\times)=\Hom_\Z(P^\gp,\bmu_N(\kappa))$, where
$\bmu_N(\kappa)\subset\kappa^\times$ is the $N$-torsion subgroup.
Hence $g_{N|U_P}$ is an \'etale Galois covering of $U_P$, whose
Galois group is naturally isomorphic to $\Hom_\Z(P^\gp,\bmu_N(\kappa))$;
therefore, $g_{N|U_P}$ is classified by a well defined continuous
representation :
\set\begin{equation}\label{eq_to-be-composed}
\pi_1(U_{P,\et},\tau)\to\Hom_\Z(P^\gp,\bmu_N(\kappa))
\end{equation}
(lemma \ref{lem_from-coh-to-fundgrp}). The corresponding
action of $\Hom_\Z(P^\gp,\bmu_N(\kappa))$ on $U_P$ can be
extracted from the construction : namely, let
$\chi:P^\gp\to\bmu_N(\kappa)$ be a given character; by
definition, the action of $\chi$ sends the geometric point
$\tau$ to the geometric point $\tau'$ whose character
$\chi_{\tau'}:P^\gp\to\kappa^\times$ is given by the rule :
$x\mapsto\chi(x)\cdot\chi_\tau$ for every $x\in P^\gp$.
Consider the automorphism $\rho_\chi$ of the $\Z$-algebra
$\Z[P^\gp]$ given by the rule $x\mapsto\chi(x)\cdot x$ for
every $x\in P^\gp$, and notice that $\rho(\tau)=\tau'$;
since the fibre functors are faithful on \'etale coverings,
we conclude that $\chi$ acts as $\Spec\,\rho_\chi$ on
$U_P$.

$\bullet$\ \
Moreover, if $\lambda:Q\to P$
is any morphism of fine and saturated monoids, clearly we have
a commutative diagram
$$
\xymatrix{
U_P \ar[rrr]^-{\Spec\,\Z[\lambda]_{|U_P}} \ar[d]_{g_{N|U_P}} & & &
U_Q \ar[d]^{g_{N|U_Q}} \\
U_P \ar[rrr]^-{\Spec\,\Z[\lambda]_{|U_P}} & & & U_Q
}$$
whence, by remark \ref{rem_from-coh-to-fungrp}(iii.a), a well
defined group homomorphism
\set\begin{equation}\label{eq_trovami}
\Hom_\Z(P^\gp,\bmu_N)\to\Hom_\Z(Q^\gp,\bmu_N)
\end{equation}
that makes commute the induced diagram :
$$
\xymatrix{
\pi_1(U_{P,\et},\xi) \ar[r] \ar[d] & \pi_1(U_Q,\xi') \ar[d] \\
\Hom_\Z(P^\gp,\bmu_N) \ar[r] & \Hom_\Z(Q^\gp,\bmu_N)
}$$
whose vertical arrows are the maps \eqref{eq_to-be-composed}.
By inspecting the constructions, it is easily seen that
\eqref{eq_trovami} is just the map $\Hom_\Z(\lambda^\gp,\bmu_N)$.

$\bullet$\ \
Next, by corollary \ref{cor_sharpy}(i) the projection
$\underline M{}_x\to P$ admits a splitting
$$
\alpha:P\to\underline M{}_x
$$
(which defines a sharp chart for $\underline M$); if $N$ is
invertible in $\cO_{\!X}$, this splitting induces a morphism of
schemes $X\to S_P$, which restricts to a morphism $U\to U_P$.
If we let $\tau$ be the image of $\xi$, we deduce a continuous
group homomorphism $\pi_1(U_\et,\xi)\to\pi_1(U_{P,\et},\tau)$
(\cite[Exp.V, \S7]{SGA1}), whose composition with
\eqref{eq_to-be-composed} yields a continuous map :
\set\begin{equation}\label{eq_yet-another-one}
\pi_1(U_\et,\xi)\to\Hom_\Z(P^\gp,\bmu_N(\kappa)).
\end{equation}
We claim that the pairing
$P^\gp\times\pi_1(U_\et,\xi)\to\bmu_N(\kappa)$ arising from
\eqref{eq_yet-another-one} agrees with \eqref{eq_first-pairing},
under the natural identification $\bmu_N(\kappa)\isom\bmu_N$.
Indeed, for given $s\in P$, let $j_s:\Z\to P^\gp$ be the map such
that $j_s(n):=s^n$ for every $n\in\Z$; by composing
\eqref{eq_yet-another-one} with $\Hom_\Z(j_s,\bmu_N(\kappa))$,
we obtain a map
$\pi_1(U_\et,\xi)\to\Hom_\Z(\Z,\bmu_N(\kappa))\isom\bmu_N(\kappa)$;
in view of the foregoing, it is easily seen that the corresponding
$\bmu_N(\kappa)$-torsor is precisely $\phi_{\alpha(s)}$, in the
notation of \eqref{eq_G_m-torsor}, whence the contention.

\begin{proposition}\label{prop_tame-quot-log}
In the situation of \eqref{subsec_first-pairing}, suppose
that $(X,\underline M)$ is a regular log scheme. Then
\eqref{eq_tame-quotient} is a surjection.
\end{proposition}
\begin{proof} From the discussion of \eqref{subsec_Second-pairing},
we are reduced to showing that the morphism
\eqref{eq_yet-another-one} is surjective, for every $N>0$ such that
$(N,p)=1$. The latter comes down to showing that the corresponding
torsor $g_{N|U_P}\times_{U_P}U:U'\to U$ is connected (remark
\ref{rem_from-coh-to-fungrp}(i)). However, the map
$\alpha:P\to\underline M{}_x$ induces a morphism of log
schemes $\psi:(X,\underline M)\to\Spec(\Z,P)$ (see
\eqref{subsec_constant-log}); we remark the following :

\begin{claim}\label{cl_Kummer-preserve-regul}
Slightly more generally, for any Kummer morphism $\nu:P\to Q$ of
monoids, with $Q$ fine, sharp and saturated, define
$(X_\nu,\underline M{}_\nu)$ as the fibre product in the cartesian
diagram:
\set\begin{equation}\label{eq_define-anew}
{\diagram
(X_\nu,\underline M{}_\nu) \ar[r] \ar[d]_{f_\nu} &
\Spec(\Z,Q) \ar[d]^{\Spec(\Z,\nu)} \\
(X,\underline M) \ar[r]^-\psi & \Spec(\Z,P)
\enddiagram}
\end{equation}
Then $(X_\nu,\underline M{}_\nu)$ is regular, and $X_\nu$
is strictly local.
\end{claim}
\begin{pfclaim} By lemma \ref{lem_prepare-for-Kummer},
$(X_\nu,\underline M{}_\nu)$ is regular at every point of the closed
fibre $X_\nu\times_X\Spec\,\kappa(x)$, and since $f_\nu$ is a
finite morphism, every point of $X_\nu$ specializes to a point of
the closed fibre, therefore $(X_\nu,\underline M{}_\nu)$ is regular
(theorem \ref{th_reg-generizes}). Moreover, $X_\nu$ is connected,
provided the same holds for the closed fibre; in turn, this follows
immediately from proposition \ref{prop_crit-Kummer}(ii). Then
$X_\nu$ is strictly local, by \cite[Ch.IV, Prop.18.5.10]{EGA4}.
\end{pfclaim}

From claim \ref{cl_Kummer-preserve-regul} we see that $X_\nu$ is
normal (corollary \ref{cor_normal-and-CM}) for every such $\nu$,
especially this holds for the $N$-Frobenius endomorphism of $P$, in
which case $U'$ is an open subset of $X_\nu$; since the latter is
connected, the same then holds for $U'$.
\end{proof}

\begin{example}\label{ex_continues-previous}
(i)\ \
In the situation of example \ref{ex_fan-andlogsch-2},
let $X\to S$ be any morphism of schemes, set
$\underline X:=X\times_S\underline S$, denote by
$(X,\underline M)$ the log scheme underlying $\underline X$,
and define $(X_{(k)},\underline M{}_{(k)})$ as the fibre
product in the cartesian diagram of log schemes :
$$
\xymatrix{
(X_{(k)},\underline M{}_{(k)}) \ar[r]^-{\pi_{(k)}} \ar[d]_{\bek_X} & 
\Spec(R,P) \ar[d]^{\Spec(R,\bek_P)} \\
(X,\underline M) \ar[r]^-\pi & \Spec(R,P)
}$$
where $\pi$ is the natural projection. Set as well
$$
\underline X{}_{(k)}:=
((X_{(k)},\underline M{}_{(k)}),T_P,\psi_P\circ\pi_{(k)}^\sharp)
$$
which is an object of $\cK_\intg$; by proposition
\ref{prop_basis-technique}(iii),  diagram \eqref{eq_yad-of-fans}
(with $T:=T_P$) underlies a commutative diagram in $\cK_\intg$ :
$$
\xymatrix{
\phi^*\underline X{}_{(k)} \ar[r] \ar[d]_{\underline g{}_X} &
\underline X{}_{(k)} \ar[d]^{(\bek_X,\bek_{T_P})} \\
\phi^*\underline X \ar[r] & \underline X
}$$
whose horizontal arrows are cartesian. Clearly this diagram is isomorphic
to $X\times_S\eqref{eq_yad-in-cK}$; we deduce that the morphism $g_X$
of log schemes underlying $\underline g{}_X$ is finite, and of Kummer type;
by proposition \ref{prop_characterize-log-etale-cov}, $g_X$ is even an 
\'etale covering, if $k$ is invertible in $\cO_{\!X}$.

(ii)\ \
Suppose furthermore, that $(X,\underline M)$ is regular;
by claim \ref{cl_Kummer-preserve-regul} it follows that
$(X_{(k)},\underline M{}_{(k)})$ is regular as well, and
then the same holds for the log schemes 
$(X_\phi,\underline M{}_\phi)$ and $(Y,\underline N)$ underlying
respectively $\phi^*\underline X$ and $\phi^*\underline X{}_{(k)}$
(proposition \ref{prop_basis-technique}(v) and theorem
\ref{th_smooth-preserve-reg}). Let now $y\in Y$ be a point of height
one in $Y$, lying in the closed subset $Y\!\setminus\!(Y,\underline N)_\tr$,
and set $x:=g_X(y)$; arguing as in example \ref{ex_Kummer-is-tamely-ram}, 
we see that $\underline N{}_y^\sharp$ and $\underline M{}_{\phi,x}^\sharp$
are both isomorphic to $\N$, and the induced map
$\cO_{\!X_\phi,x}\to\cO_{Y,y}$ is an extension of discrete valuation
rings, whose corresponding extension of valued fields is finite.
Denote by $\Gamma_x\to\Gamma_y$ the associated extension of value
groups; in view of  \eqref{eq_finnegan} we see that the ramification
index $(\Gamma_y:\Gamma_x)$ equals $k$.
\end{example}

\sset\subsubsection{}\label{subsec_ex-tame-covers}
Resume the situation of \eqref{subsec_Second-pairing}.
Every (finite) discrete quotient map
$$
\rho:
\underline M{}_x^{\gp\vee}\otimes_\Z\prod_{\ell\neq p}\Z_\ell(1)
\to G
$$
corresponds, via composition with \eqref{eq_tame-quotient}, to
a $G_U$-torsor, which can be explicitly constructed as follows.
Set $P:=\underline M{}_x^\sharp$, pick a splitting $\alpha$ as
in \eqref{subsec_Second-pairing}, choose an integer $N>0$ large
enough, so that $(N,p)=1$, and $\rho$ factors through a group
homomorphism
$$
\bar\rho:\Hom_\Z(P^\gp,\bmu_N(\kappa))\to G
$$
and set $H:=\Ker\,\bar\rho$. Now, via \eqref{eq_tame-quotient},
the quotient $G':=\Hom_\Z(P^\gp,\bmu_N(\kappa))$ corresponds to
the $G'_U$-torsor on $U$ obtained from $g_N:U_P\to U_P$ via
pull-back along the morphism $U\to U_P$ given by the chart
$\alpha$ (notation of \eqref{subsec_Second-pairing}). According
to remark \ref{rem_from-coh-to-fungrp}(ii), the sought $G_U$
is therefore isomorphic to the one obtained from the quotient
$\bar g_N:U_P/H\to U_P$, via pull-back along the same morphism.
To exhibit such quotient, consider the exact sequence:
$$
0\to\Hom_\Z(G,\kappa^\times)\to
P^\gp\otimes_\Z\Z/N\Z\to\Hom_\Z(H,\kappa^\times)\to 0.
$$
Define $Q^\gp\subset P^\gp$ as the kernel of the induced map
$P^\gp\to\Hom_\Z(H,\kappa^\times)$, and set $Q:=Q^\gp\cap P$. By
construction, the $N$-Frobenius endomorphism of $P$ factors through
an injective map $\nu:P\to Q$ and the inclusion map $j:Q\to P$,
so $g_N$ factors as a composition :
$$
U_P\to U_Q\xrightarrow{\ h\ } U_P.
$$
The maps on geometric points $U_P(\kappa)\to U_Q(\kappa)\to U_P(\kappa)$
induced by $g_N$ and $h$ correspond to $j^{\gp*}$ and respectively
$\nu^{\gp*}$ in the resulting commutative diagram
$$
\xymatrix{
& & & 0 \ar[d] \\
0 \ar[r] & H \ar[r] \ddouble &
\Hom_\Z(P^\gp,\bmu_N(\kappa)) \ar[r] \ar[d] & G \ar[d] \ar[r] & 0 \\
0 \ar[r] & H \ar[r] & \Hom_\Z(P^\gp,\kappa^\times) \ar[r]^-{j^{\gp*}} &
\Hom_\Z(Q^\gp,\kappa^\times) \ar[d]^-{\nu^{\gp*}} \ar[r] & 0 \\
& & & \Hom_\Z(P^\gp,\kappa^\times) \ar[d] \\
& & & 0
}$$
whose rows and column are exact.
Hence, let $\tau$ and $\chi_\tau:P^\gp\to\kappa^\times$ be as in
\eqref{subsec_Second-pairing}; the fibre $h^{-1}(\tau)$ corresponds
to the set of all characters $\chi_{\tau'}:Q^\gp\to\kappa^\times$
whose restriction $\chi_{\tau'}\circ\nu^\gp$ to $P^\gp$ agrees with
$\chi_\tau$. Then, with the notation of \eqref{eq_define-anew},
we conclude that the restriction
$(f_\nu)_\tr:(X_\nu,\underline M_\nu)_\tr\to U$ is the sought
$G_U$-torsor. Notice as well that $f_\nu$ is an \'etale covering
of $(X,\underline M)$.

$\bullet$\ \
Here is another handier description of the same submonoid
$Q$ of $P$. Notice that $\rho$ is the same as a group
homomorphism
$$
\rho^\dagger:P^{\gp\vee}\to\Hom_\Z(\textstyle{\prod_{\ell\neq p}}\Z_\ell(1),G)
$$
and let $L:=\Ker\,\rho^\dagger$. Fix a generator $\zeta_N$ of
$\bmu_N(\kappa)$; by definition
$$
\begin{aligned}
Q^\gp=\: & \{x\in P^\gp~|~\text{$t(x)\otimes\xi=0$ for all
$t\in P^{\gp\vee}$ and $\xi\in\bmu_N$ such that
$\bar\rho(t\otimes\xi)=0$}\} \\
=\: & \{x\in P^\gp~|~\text{$t(x)\otimes\zeta_N=0$ for all
$t\in P^{\gp\vee}$ such that $\bar\rho(t\otimes\zeta_N)=0$}\} \\
=\: & \{x\in P^\gp~|~\text{$t(x)\in N\Z$ for all
$t\in L$}\}.
\end{aligned}
$$
On the other hand, notice that the $N$-Frobenius of $P^{\gp\vee}$
factors through a map $\beta:P^{\gp\vee}\to L$ and the inclusion
map $i:L\to P^{\gp\vee}$, and $\beta\circ i$ is the $N$-Frobenius
map of $L$. Let $\omega:P^\gp\isom(P^{\gp\vee})^\vee$ be the
natural isomorphism. We may then write
$$
\begin{aligned}
Q^\gp=\: & \{x\in P^\gp~|~
\text{$\omega(x)(t)\in N\Z$ for all $t\in L$}\} \\
=\: & \{x\in P^\gp~|~
\text{$\omega(x)\circ i\in N\cdot L^{\vee}=i^\vee\circ\beta^\vee(L^\vee)$}\} \\
=\: & \omega^{-1}(\Img\,\beta^\vee).
\end{aligned}
$$
In other words, $\omega$ induces a natural isomorphism
$Q^\gp\isom L^\vee$, that identifies $\beta^\vee$ with
the map $j^\gp:Q^\gp\to P^\gp$, and then necessarily also
the map $i^\vee$ with $\nu$. Now, the image of $Q$ inside
$L^\vee$ can be recovered just as $L^\vee\cap P^{\vee\vee}_\Q$
(the intersection here takes place in $L^\vee_\Q$, which
contains $(P^{\gp\vee})^\vee$, via the injective map
$i^\vee$ : details left to the reader).

\sset\subsubsection{}\label{subsec_supply}
Our chief supply of tamely ramified coverings comes from the
following source. Let $f:(Y,\underline N)\to(X,\underline M)$
be an \'etale morphism of log schemes, whose underlying morphism
of schemes is finite, and suppose that $(X,\underline M)$ is regular.
Then $(Y,\underline N)$ is regular, and both $X$ and $Y$ are normal
schemes (theorem \ref{th_smooth-preserve-reg} and corollary
\ref{cor_normal-and-CM}).
Moreover, lemma \ref{lem_Kummer-fans}(ii) and proposition
\ref{prop_characterize-log-etale-cov} imply that $(Y,\underline
N)_\tr=f^{-1}(X,\underline M)_\tr$, hence the restriction of $f$
$$
f_\tr:(Y,\underline N)_\tr\to(X,\underline M)_\tr
$$
is a finite \'etale morphism of schemes (corollary
\ref{cor_undercover}(i)). Furthermore, it follows easily from
corollary \ref{cor_same-height}(i) that the closed subset
$X\!\setminus\!(X,\underline M)_\tr$ is a union of irreducible
closed subsets of codimension $1$ in $X$ (and the union is locally
finite on the Zariski topology of $X$); the same holds also for
$Y\!\setminus\!(Y,\underline N)_\tr$, especially $Y$ is the
normalization of $(Y,\underline N)_\tr$ over $X$.
Finally, example \ref{ex_Kummer-is-tamely-ram} implies that $f_\tr$
is tamely ramified along $X\!\setminus\!(X,\underline M)_\tr$. It is
then clear that the rule $f\mapsto f_\tr$ defines a functor :
$$
F_{(X,\underline M)}:\bCov(X,\underline
M)\to\bTame(X,(X,\underline M)_\tr).
$$
It follows easily from remark \ref{rem_fully-regular} that
$F_{(X,\underline M)}$ is fully faithful; therefore, any choice
of a geometric point $\xi$ of $(X,\underline M)_\tr$ determines
a surjective group homomorphism :
\set\begin{equation}\label{eq_cuneo}
\pi_1((X,\underline M)_{\tr,\et},\xi)\to
\pi_1((X,\underline M)_\et,\xi)
\end{equation}
(\cite[Exp.V, Prop.6.9]{SGA1}).

\begin{proposition}\label{prop_fund-cuneo}
In the situation of \eqref{subsec_first-pairing}, suppose that
$(X,\underline M)$ is a regular log scheme. Then the map
\eqref{eq_tame-quotient} factors through \eqref{eq_cuneo},
and induces an isomorphism :
\set\begin{equation}\label{eq_madonnadellolmo}
\pi_1((X,\underline M)_\et,\xi)\isom
\underline M{}_{\bar x}^{\gp\vee}\otimes_\Z\prod_{\ell\neq p}\Z_\ell(1).
\end{equation}
\end{proposition}
\begin{proof} The discussion of \eqref{subsec_ex-tame-covers}
shows that \eqref{eq_tame-quotient} factors through \eqref{eq_cuneo},
and proposition \ref{prop_tame-quot-log} implies that
\eqref{eq_madonnadellolmo} is surjective. Next, let
$f:(Y,\underline N)\to(X,\underline M)$ be an \'etale covering;
according to proposition \ref{prop_characterize-log-etale-cov},
$f$ admits a Kummer chart $(\omega_P,\omega_Q,\theta)$ with
$Q$ fine, sharp and saturated, such that the order $k$ of
$\Coker\,\theta^\gp$ is invertible in $\cO_Y$. Moreover,
the induced morphism of $X$-schemes
$Y\to X\times_{\Spec\,\Z[P]}\Spec\,\Z[Q]$ is an isomorphism.
With this notation, denote by $G$ the cokernel of the induced
group homomorphism $\Hom_\Z(\theta^\gp,\bek_{P^\gp})$; there
follows a map
$\rho:P^\gp\otimes_\Z\prod_{\ell\neq p}\Z_\ell(1)\to G$,
and by inspecting the construction, one may check that the
corresponding $G_U$-torsor, constructed as in
\eqref{subsec_ex-tame-covers}, is isomorphic to $f\times_X\one_U$
(details left to the reader). This shows the injectivity of
\eqref{eq_madonnadellolmo}, and completes the proof of the
proposition.
\end{proof}

The following result is the logarithmic version of the
classical Abhyankar's lemma.

\begin{theorem}\label{th_Abhyankar_log}
The functor $F_{(X,\underline M)}$ is an equivalence.
\end{theorem}
\begin{proof} First we show how to reduce to the case where
$\tau=\et$.

\begin{claim}\label{cl_redux-to-etale-F}
Let $(X_\Zar,\underline M)$ be a regular log structure on the
Zariski site of $X$, and suppose that the theorem holds for
$\tilde u{}^*(X_\Zar,\underline M)$. Then the theorem holds for
$(X_\Zar,\underline M)$ as well.
\end{claim}
\begin{pfclaim}
Let $g:V\to(X_\Zar,\underline M)_\tr$ be an \'etale covering whose
normalization over $X$ is tamely ramified along
$X\!\setminus\!(X_\Zar,\underline M)_\tr$. By assumption, there
exists an \'etale covering $f:(Y_\et,\underline N)\to\tilde
u{}^*(X_\Zar,\underline M)$ such that $f_\tr=g$; by lemma
\ref{lem_Kummer-descends}, we may then find a morphism of log
schemes of Kummer type $f_\Zar:\tilde u_*(Y_\et,\underline
N)\to(X_\Zar,\underline M)$ such that $\tilde u{}^*f_\Zar=f$,
therefore $f_\Zar$ is an \'etale covering of $(X_\Zar,\underline M)$
(corollary \ref{cor_undercover}(iii)), and clearly $(f_\Zar)_\tr=g$,
so the functor  $F_{(X_\Zar,\underline M)}$ is essentially
surjective. Full faithfulness for the same functor is derived
formally from the full faithfulness of the functor $F_{\tilde
u{}^*(X_\Zar,\underline M)}$, and that of the functor
\eqref{eq_pulbak-kummer} : details left to the reader.
\end{pfclaim}

Henceforth we assume that $\tau=\et$.

\begin{claim}\label{cl_first-redux}
Let $g:V\to(X,\underline M)_\tr$ be an object of
$\bTame(X,(X,\underline M)_\tr)$. Then $g$ lies in the
essential image of $F_{(X,\underline M)}$ if (and only if) there
exists an \'etale covering $(U_\lambda\to X~|~\lambda\in\Lambda)$ of
$X$, such that, for every $\lambda\in\Lambda$, the \'etale covering
$g\times_XU_\lambda$ lies in the essential image of
$F_{(U_\lambda,\underline M_{|U_\lambda})}$.
\end{claim}
\begin{pfclaim} We have to exhibit a finite \'etale covering
$f:(Y,\underline N)\to(X,\underline M)$ such that $f_\tr=g$.
However, given such $f$, theorem \ref{th_smooth-preserve-reg} and
corollary \ref{cor_normal-and-CM} imply that $Y$ is the
normalization of $V$ over $X$, and proposition
\ref{prop_log-struct-is-fixed} says that $\underline
N=j_*\cO_V^\times\cap\cO_Y$, where $j:V\to Y$ is the open immersion;
then $\log f:f^*\underline M\to\underline N$ is completely
determined by $f$, hence by $g$. Thus, we come down to showing that
:
\begin{enumerate}
\item
$\underline N:=j_*\cO_V^\times\cap\cO_Y$ is a regular log structure
on the normalization $Y$ of $V$ over $X$.
\item
The unique morphism $f:(Y,\underline N)\to(X,\underline M)$ is an
\'etale covering.
\end{enumerate}
However, \cite[\S33, Lemma 1]{Mat} implies that $Y$ is finite over
$X$. To show (ii), it then suffices to prove that each restriction
$f_\lambda:=f\times_X\one_{U_\lambda}$ is \'etale (proposition
\ref{prop_sorite-smooth}(iii)). Likewise, (i) holds, provided the
restriction $\underline N_{|Y_\lambda}$ is a regular log structure
on $Y_\lambda:=Y\times_XU_\lambda$, for every $\lambda\in\Lambda$.
Let $j_\lambda:V_\lambda:=V\times_XU_\lambda\to Y_\lambda$ be the
induced open immersion. It is clear that $\underline
N_{|Y_\lambda}=j_{\lambda*}\cO_{V_\lambda}^\times\cap\cO_{Y_\lambda}$,
and $Y_\lambda$ is the normalization of $V_\lambda$ over $U_\lambda$
(\cite[Ch.IV, Prop.17.5.7]{EGA4}); moreover, $(U_\lambda,\underline
M_{|U_\lambda})$ is again regular, so our assumption implies that
$(Y_\lambda,\underline N_{|Y_\lambda})$ is the unique regular log
structure on $Y_\lambda$ whose trivial locus is $V_\lambda$, and
that $f_\lambda$ is indeed \'etale.
\end{pfclaim}

\begin{claim}\label{cl_reduce-to-strict-hens}
If $F_{(X(\xi),\underline M(\xi))}$ is essentially surjective for
every geometric point $\xi$ of $X$, the same holds for
$F_{(X,\underline M)}$.
\end{claim}
\begin{pfclaim} Indeed, let $g:V\to(X,\underline M)_\tr$ be an object of
$\bTame(X,(X,\underline M)_\tr)$, and $\xi$ any geometric
point of $X$. By claim \ref{cl_first-redux} it suffices to find an
\'etale neighborhood $U\to X$ of $\xi$, such that $g\times_X\one_U$
lies in the essential image of $F_{(U,\underline M_{|U})}$. Denote
by $Y$ the normalization of $X$ in $V$, which is a finite
$X$-scheme (\cite[\S33, Lemma 1]{Mat}). The scheme
$Y(\xi):=Y\times_XX(\xi)$ decomposes as a finite disjoint union of
strictly local open and closed subschemes $Y_1(\xi),\dots,Y_n(\xi)$,
and we may find an \'etale neighborhood $U\to X$ of $\xi$, and a
decomposition of $Y\times_XU$ by open and closed subschemes
$Y_1,\dots,Y_n$, with isomorphisms of $X(\xi)$-schemes
$Y_i(\xi)\isom Y_i\times_XX(\xi)$ for every $i=1,\dots,n$
(\cite[Ch.IV, Cor.8.3.12]{EGAIV-3}). We are then reduced to showing
that all the restrictions $Y_i\to U$ lie in the essential image of
$F_{(U,\underline M_{|U})}$. Thus, we may replace $X$ by $U$, and
$Y$ by any $Y_i$, after which we may assume that $Y(\xi)$ is
strictly local. Proposition \ref{prop_second-crit} and theorem
\ref{th_reg-generizes} imply that $(X(\xi),\underline M(\xi))$ is a
regular log scheme. Clearly $g_\xi:=g\times_X\one_{X(\xi)}$ is an
object of $\bTame(X(\xi),(X(\xi),\underline M(\xi))_\tr)$, so
by assumption there exists a finite \'etale covering
$h:(Z,\underline N)\to(X(\xi),\underline M(\xi))$ such that
$h_\tr=g_\xi$. By theorem \ref{th_smooth-preserve-reg} and corollary
\ref{cor_normal-and-CM} we know that $Z$ is a normal scheme, and we
deduce that $Z=Y(\xi)$ (\cite[Ch.IV, Prop.17.5.7]{EGA4}).

According to proposition \ref{prop_characterize-log-etale-cov}, the
morphism $h$ admits a fine and saturated Kummer chart
$(\omega_P,\omega_Q,\theta:P\to Q)$, such that the order $d$ of
$\Coker\,\theta^\gp$ is invertible in $\cO_{\!X,\xi}$, and such that
the induced map $Y(\xi)\to X(\xi)\times_{\Spec\,P}\Spec\,Q$ is an
isomorphism. By proposition \ref{prop_2-colim-for-logs}, there
exist an \'etale neighborhood $U\to X$ of $\xi$, and a coherent log
structure $\underline N'$ on $Y':=Y\times_XU$ with an isomorphism
$Y(\xi)\times_{Y'}(Y',\underline N')\isom(Y(\xi),\underline N)$.
Moreover, let $h':Y'\to U$ be the projection; after shrinking $U$,
the map $\log h$ descends to a morphism of log structures
$h^{\prime*}\underline M_{|U}\to\underline N'$, whence a morphism
$(h',\log h'):(Y',\underline N')\to(U,\underline M_{|U})$ of log
schemes, such that $h'\times_U\one_{X(\xi)}=h$. After further
shrinking $U$, we may also assume that $h'$ admits a Kummer chart
$(\omega'_P,\omega'_Q,\theta)$ (corollary
\ref{cor_descend-chart-from-infty}), that $d$ is invertible in
$\cO_{\!U}$, and that the induced morphism $Y'\to
U\times_{\Spec\,P}\Spec\,Q$ is an isomorphism (\cite[Ch.IV,
Cor.8.8.24]{EGAIV-3}). Then $h'$ is an \'etale covering, by
proposition \ref{prop_characterize-log-etale-cov}, and by
construction, $F_{(U,\underline M_{|U})}(h')=g\times_X\one_U$, as
desired.
\end{pfclaim}

\begin{claim}\label{cl_replace-by-g_phi}
Assume that $X$ is strictly local, denote by $x$ the
closed point of $X$, set $U:=(X,\underline M)_\tr$,
and $P:=\underline M{}^\sharp_x$. Let $h:V\to U$ be any connected
non-empty \'etale covering, tamely ramified along $X\!\setminus\!U$;
then $h$ lies in the essential image of $F_{(X,\underline M)}$, provided
there exist a fine, sharp and saturated monoid $Q$, and a morphism
$\nu:P\to Q$ of Kummer type, such that
$$
h_\nu:=h\times_XX_\nu:V\times_XX_\nu\to U\times_XX_\nu
$$
admits a section (notation of \eqref{eq_define-anew}).
\end{claim}
\begin{pfclaim} Suppose first that the order of $\Coker\,\nu^\gp$
is invertible in $\cO_{\!X}$; in this case, $f_\nu$ is an \'etale
covering of log schemes (proposition
\ref{prop_characterize-log-etale-cov}), hence $(f_\nu)_\tr$ is an
\'etale covering of the scheme $U$. Then, by composing a section of
$h_\nu$ with the projection $V\times_XX_\nu\to V$ we deduce a
morphism $U\times_XX_\nu\to V$ of \'etale coverings of $U$. Such
morphism shall be open and closed, hence surjective, since $V$ is
connected; hence the $\pi_1(U,\xi)$-set $h^{-1}(\xi)$ will be a
quotient of $f^{-1}_\nu(\xi)$, on which $\pi_1(U,\xi)$ acts through
\eqref{eq_tame-quotient}, as required. Next, for a general morphism
$\nu$ of Kummer type, let $L\subset Q^\gp$ be the largest subgroup
such that $\nu^\gp(P^\gp)\subset L$, and $(L:P^\gp)$ is invertible
in $\cO_{\!X}$. Set $Q':=L\cap Q$, and notice that
$Q^{\prime\gp}=L$. Indeed, every element of $L$ can be written in
the form $x=b^{-1}a$, for some $a,b\in Q$; then choose $n>0$ such
that $b^n\in\nu P$, write $x=b^{-n}\cdot(b^{n-1}a)$ and remark that
$b^{-n},b^{n-1}a\in Q'$. The morphism $\nu$ factors as the
composition of $\nu':P\to Q'$ and $\psi:Q'\to Q$, and therefore
$f_\nu$ factors through a morphism $f_\psi:X_\nu\to X_{\nu'}$. In
view of the previous case, we are reduced to showing that the
morphism $h_{\nu'}$ already admits a section, hence we may replace
$(X,\underline M)$ by $(X_{\nu'},\underline M{}_{\nu'})$ (which is
still regular and strictly local, by claim
\ref{cl_Kummer-preserve-regul}), $h$ by $h_{\nu'}$, and $\nu$ by
$\psi$, after which we may assume that the order of
$\Coker\,\nu^\gp$ is $p^m$ for some integer $m>0$, where $p$ is the
characteristic of the residue field $\kappa(x)$, and then we need to
show that $h$ already admits a section.

Next, using again claim \ref{cl_Kummer-preserve-regul} and an easy
induction, we may likewise reduce to the case where $m=1$. Say that
$X=\Spec\,A$, and suppose first that $A$ is a $\F_p$-algebra (where
$\F_p$ is the finite field with $p$ elements), so that
$X_\nu=\Spec\,A\otimes_{\F_p[P]}\F_p[Q]$; it is easily seen that
the ring homomorphism $\F_p[\nu]:\F_p[P]\to\F_p[Q]$ is invertible
up to $\Phi$, in the sense of \cite[Def.3.5.8(i)]{Ga-Ra}.
Especially, $\Spec\,\F_p[\nu]$ is integral, surjective and
radicial, hence the same holds for $f_\nu$, and therefore the
morphism of sites :
$$
f_\nu^*:X_\et\to X_{\nu,\et}
$$
is an equivalence of categories (lemma \ref{lem_replace}(i)). It
follows that in this case, $h$ admits a section if and only if the
same holds for $h_\nu$.

Next, suppose that the field of fractions $K$ of $A$ has
characteristic zero; we may write $X_\nu\times_X\Spec\,K=\Spec\,K_\nu$
and $V\times_X\Spec\,K=\Spec\,E$ for two field extensions $K_\nu$
and $E$ of $K$, such that $[K_\nu:K]=p$, and the section of
$h_\nu$ yields a map $E\to K_\nu$ of $K$-algebras. Therefore we
have either $E=K$ (in which case $V=U$, and then we are done), or
else $E=K_\nu$, in which case $V=(X_\nu,\underline M{}_\nu)_\tr$,
since both these schemes are normal and finite over $U$.

Hence, let us assume that $V=(X_\nu,\underline M{}_\nu)$, and pick
any point $\eta\in X$ of codimension one, whose residue field
$\kappa(\eta)$ has characteristic $p$. Then
$X_\nu\times_XX(\eta)=\Spec\,B$, where
$B:=\cO_{\!X,\eta}\otimes_{\Z[P]}\Z[Q]$, and
$\bar B:=B\otimes_A\kappa(\eta)=\kappa(\eta)\otimes_{\F_p[P]}\F_p[Q]$.
Since the map $\F_p[P]\to\F_p[Q]$ is invertible up to $\Phi$, we
easily deduce that $\bar B$ is local, and its residue field is a
purely inseparable extension of $\kappa(\eta)$, say of degree $d$.
Therefore $f_\nu^{-1}\eta$ consists of a single point $\eta'$, and
$\cO_{\!X_\nu,\eta'}\simeq B$ is a normal and finite
$\cO_{\!X,\eta}$-algebra, hence it is a discrete valuation ring.
Let $e$ denote the ramification index of the extension
$\cO_{\!X,\eta}\to\cO_{\!X_\nu,\eta'}$; then $ed=p$
(\cite[Ch.VI, \S8, n.5, Cor.1]{BouAC}). Suppose that $e=p$; in
view of \cite[Lemma 6.2.5]{Ga-Ra}, the ramification index of the
induced extension of strict henselizations
$\cO_{\!X,\eta}^\sh\to\cO_{\!X_\nu,\eta'}^\sh$ is still equal to
$p$, which is impossible, since $h$ is tamely ramified along
$X\!\setminus\!U$. In case $e=1$, we must have $d=p$, and since the
residue field $\kappa(\eta)^\mathrm{s}$ of $\cO_{\!X,\eta}^\sh$ is a
separable closure of $\kappa(\eta)$, it follows easily that the
residue field of $\cO_{\!X_\nu,\eta'}^\sh$ must be a purely
inseparable extension of $\kappa(\eta)^\mathrm{s}$ of degree $p$,
which again contradicts the tameness of $h$.
\end{pfclaim}

\begin{claim}\label{cl_up-to-cod-2}
If $(X,\underline M)=(X,\underline M)_2$, then $F_{(X,\underline
M)}$ is essentially surjective (notation of definition
\ref{def_trivial-locus}(i)).
\end{claim}
\begin{pfclaim} By claim \ref{cl_reduce-to-strict-hens}, we may assume
that $X$ is strictly henselian. In this case, let $U$, $P$ and
$h:V\to U$ be as in claim \ref{cl_replace-by-g_phi}. By assumption
$d:=\dim P\leq 2$, and we have to show that $h\times_XX_\nu$ admits
a section, for a suitable Kummer morphism $\nu:P\to Q$ (notation
of \eqref{eq_define-anew}). If $d=0$, then $P=\{1\}$, in which
case $U=X$ is strictly local, so its fundamental group is trivial,
and there is nothing to prove. In case $d=1$, then $P\simeq\N$
(theorem \ref{th_structure-of-satu}(ii)), in which case $X$ is a
regular scheme (corollary \ref{cor_more-precisely}), and
$X\!\setminus\!U$ is a regular divisor (remark
\ref{rem_simplex-is-normal-cross}) and then the assertion follows
from the classical Abhyankar's lemma (\cite[Exp.XIII,
Prop.5.2]{SGA1}). For $d=2$, we may find $e_1,e_2\in P$, and an
integer $N>0$ such that
$$
\N e_1\oplus\N e_2\subset P\subset
P':=\N\frac{e_1}{N}\oplus\N\frac{e_2}{N}
$$
(see example \ref{ex_satu-dim-two}(i)). Especially, the inclusion
$\nu:P\to P'$ is a morphism of Kummer type. Since a composition of
morphisms of Kummer type is obviously of Kummer type, claims
\ref{cl_replace-by-g_phi} and \ref{cl_Kummer-preserve-regul} imply
that we may replace $(X,\underline M)$ by $(X_\nu,\underline M{}_\nu)$
and $h$ by $h_\nu$ (notation of \eqref{eq_define-anew}), after which
we may assume that $P$ is isomorphic to $\N^{\oplus 2}$.
In this case, $X$ is again regular and $X\!\setminus\!U$ is a
strict normal crossings divisor, so the assertion follows again
from \cite[Exp.XIII, Prop.5.2]{SGA1}.
\end{pfclaim}

\begin{claim}\label{cl_tame-to-tame}
In the situation of \eqref{subsec_subdiv-and-pi_one}, suppose
that $(Y,\underline N)$ is regular. Then the functor
$f^*_\tr:\bCov(Y,\underline N)_\tr\to\bCov(Y',\underline N')_\tr$
restricts to an equivalence :
\set\begin{equation}\label{eq_equivalence-Maria}
\bTame(Y,(Y,\underline N)_\tr)\isom
\bTame(Y',(Y',\underline N')_\tr).
\end{equation}
\end{claim}
\begin{pfclaim} Arguing as in the proof of proposition
\ref{prop_from-proper-base-change}, we may reduce to the
case where $\phi:T'\to T$ is the saturated blow up of an
ideal generated by two elements of $\Gamma(T,\cO_T)$.

Notice that $(f,\log f)$ is an \'etale morphism (proposition
\ref{prop_basis-technique}(v)), and it restricts to an isomorphism
$f_\tr:(Y',\underline N')_\tr\isom(Y,\underline N)_\tr$ (remark
\ref{rem_basic-technique}). Especially, $(Y',\underline N')$ is
regular, and $f^*_\tr$ is trivially an equivalence from the \'etale
coverings of $(Y,\underline N)_\tr$ to those of $(Y',\underline N')_\tr$.

Let $g:V\to(Y,\underline N)_\tr$ be an object of
$\bTame(Y,(Y,\underline N)_\tr)$. By claim \ref{cl_up-to-cod-2}
and lemma \ref{lem_reduce-to-et}, the functor
$F_{\tilde u{}^*(Y,\underline N)_2}$ is essentially surjective,
and then the same holds for $F_{(Y,\underline N)_2}$, in view of
claim \ref{cl_redux-to-etale-F}; hence we may find an \'etale
covering $(W,\underline Q)\to(Y,\underline N)_2$, and an
isomorphism $V\isom(Y,\underline N)_\tr\times_YW$ of
$(Y,\underline N)_\tr$-schemes. On the other hand, example
\ref{ex_blow-ups}(ii) shows that $\phi$ restricts to a morphism
$T'_1\to T_2$ (notation of \eqref{subsec_height-in-T}), consequently
$f$ restricts to a morphism $(Y',\underline N')_1\to(Y,\underline N)_2$.
There follows a well defined \'etale covering :
$$
(Y',\underline N')_1\times_{(Y,\underline N)_2}(W,\underline Q)
\to(Y',\underline N')_1\times_{Y'}(Y',\underline N').
$$
whose image under $F_{(Y',\underline N')}$ is $f^*_\tr(g)$.
However, $(Y',\underline N')_1$ contains all the points of
height one of $Y'\!\setminus\!(Y',\underline N')_\tr$
(corollary \ref{cor_same-height}(i)), hence $f^*_\tr(g)$ is tamely
ramified along $Y'\!\setminus\!(Y',\underline N')_\tr$ (see
\eqref{subsec_supply}). This shows that \eqref{eq_equivalence-Maria}
is well defined, and clearly this functor is fully faithful, since
the same holds for $f^*_\tr$.
\end{pfclaim}

Suppose again that $(X,\underline M)$ is strictly henselian, define
$P$, $U$ and $h:V\to U$ as in claim \ref{cl_replace-by-g_phi}, and
set $T_P:=(\Spec\,P)^\sharp$; it is easily seen that the counit of
adjunction
$$
\tilde u{}^*\tilde u_*(X,\underline M)\to(X,\underline M)
$$
is an isomorphism (cp. \eqref{subsec_down-from-log}), hence
$\tilde u_*(X,\underline M)$ is regular (lemma \ref{lem_reduce-to-et}).
Let
$$
\underline X:=(\tilde u_*(X,\underline
M),T_P,\pi_X)
$$
be the object of $\cK$ arising from $\tilde u_*(X,\underline M)$,
as in \eqref{subsec_fan-of-a-logsch}; by theorem
\ref{th_resolution-Zar} (and its proof), we may find an
integral proper simplicial subdivision $\phi:F\to T_P$,
such that the morphism of schemes underlying
$(f,\phi):\phi^*\underline X\to\underline X$ is a resolution of
singularities for $X$.
Let $k$ be the ramification index of the covering $h$; we define
$(X_{(k)},\underline M{}_{(k)})$, $(X_\phi,\underline M{}_\phi)$, 
$(Y,\underline N)$, $g_X$ and $\bek_X$ as in example
\ref{ex_continues-previous} : this makes sense, since, in the
current situation, $\underline X$ is isomorphic to
$X\times_S\underline S$ (where $\underline S$ is defined as in
example \ref{ex_fan-andlogscheme}(i)). Moreover, by the same token,
the morphism of schemes $Y\to X_{(k)}$ is a resolution of
singularities.

Set as well $U_\phi:=(X_\phi,\underline M{}_\phi)_\tr$ and
$U_{(k)}:=(X_{(k)},\underline M{}_{(k)})_\tr$; by claim
\ref{cl_tame-to-tame} and lemma \ref{lem_tameness-pullsback}(i),
we have an essentially commutative diagram of functors :
$$
\xymatrix{
\bTame(X,U) \ar[r]^-{f^*_\tr} \ar[d]_{\bek_X^*} & 
\bTame(X_\phi,U_\phi) \ar[d]^{g_X^*} \\
\bTame(X_{(k)},U_{(k)}) \ar[r] & \bTame(Y,(Y,\underline N)_\tr)
}$$
whose horizontal arrows are equivalences. In light of claim
\ref{cl_replace-by-g_phi}, we are reduced to showing that
$\bek_X^*(h)$ admits a section. Set $h_\phi:=f^*_\tr(h)$;
then it suffices to show that $g_X^*(h_\phi)$ admits a section,
and notice that the ramification index of $h_\phi$ along
$X_\phi\!\setminus\!U_\phi$ divides $k$ (lemma
\ref{lem_tameness-pullsback}(ii)).

Let $\bar V\to X_\phi$ be the normalization of $h_\phi$ over
$X_\phi$, and $w$ a geometric point of $\bar V\times_{X_\phi}Y$
whose support has height one; denote by $y$ (resp. $v$, resp. $x$)
the image of $w$ in $Y$ (resp. in $\bar V$, resp. in $X_\phi$),
and suppose that the support of $x$ lies in $U_\phi$. Denote
by $K^\sh(x)$ the field of fractions of the strict henselization
$\cO^\sh_{\!X,x}$ of $\cO_{\!X,x}$, and define likewise $K^\sh(w)$,
$K^\sh(y)$ and $K^\sh(v)$. There follows a commutative diagram
of inclusions of valued fields :
$$
\xymatrix{ K^\sh(x) \ar[r] \ar[d] & K^\sh(y) \ar[d] \\
                  K^\sh(v) \ar[r] & K^\sh(w)
}$$
which is cocartesian in the category of fields ({\em i.e.}
$K^\sh(w)$ is the compositum of $K^\sh(y)$ and $K^\sh(v)$ : cp.
the proof of lemma \ref{lem_tameness-pullsback}(i)).
By the foregoing, the ramification index of $\cO^\sh_{\bar V,v}$
over $\cO^\sh_{\!X_\phi,x}$ divides $k$; on the other hand, example
\ref{ex_continues-previous}(ii) shows that the ramification index
of $\cO^\sh_{Y,y}$ over $\cO^\sh_{\!X_\phi,x}$ equals $k$.
It then follows ({\em e.g.} from \cite[Claim 6.2.15]{Ga-Ra}) that
$K^\sh(v)\subset K^\sh(y)$, and therefore $K^\sh(w)=K^\sh(y)$;
this shows that the ramification index of $g_X^*(h_\phi)$ along
$Y\!\setminus\!(Y,\underline N)_\tr$ equals $1$, therefore $g_X^*(h_\phi)$
is the restriction of an \'etale covering $\bar h$ of $Y$ (lemma
\ref{lem_tameness-pullsback}(iii)). By proposition
\ref{prop_from-proper-base-change}(ii), $\bar h$ is the pull-back
of an \'etale covering of $X_{(k)}$. Since $X_{(k)}$ is strictly local
(claim \ref{cl_Kummer-preserve-regul}), it follows that $\bar h$
admits a section, hence the same holds for $g_X^*(h_\phi)$, as
required.
\end{proof}

\begin{remark} A proof of theorem  \ref{th_Abhyankar_log}
similar to the one given here can be found in \cite[\S2.3]{Mo}.
\end{remark}

\sset\subsubsection{}\label{subsec_punctured}
As an application of theorem \ref{th_Abhyankar_log}, we shall
determine the fundamental group of the ``punctured'' scheme
obtained by removing the closed point from a strictly local
regular log scheme $(X,\underline M)$ (with non-trivial log
structure) of dimension $\geq 2$. Indeed, let $x\in X$ be the
closed point, and set $r:=\dim\underline M{}_x$. Also, let
$\bar y$ be any geometric point of $X$, localized at a point
$y\in U:=X\setminus\{x\}$, and $\xi$ a geometric point of $X$
localized at the maximal point. Let $p$ (resp. $p'$) denote the
characteristic exponent of $\kappa(x)$ (resp. of $\kappa(y)$).
Pick any lifting of $\xi$ to a geometric point $\xi_y$ of
$X(\bar y)$; since $(X(\bar y),\underline M(\bar y))$ is
regular (theorem \ref{th_reg-generizes}), according to proposition
\ref{prop_fund-cuneo}, the vertical arrows of the commutative
diagram in \eqref{subsec_funct-standard-shit} factor through
the surjections \eqref{eq_cuneo}, and we get a commutative
diagram of group homomorphisms :
$$
\xymatrix{ \pi_1((X(\bar y),\underline M(\bar y))_\et,\xi_y)
\ar[r]^-{\phi_y} \ar[d] & \pi_1((X,\underline M)_\et,\xi) \ar[d] \\
\underline M{}^{\gp\vee}_{\bar y}\otimes_\Z\prod_{\ell\neq p'}\Z_\ell(1)
\ar[r] &
\underline M{}^{\gp\vee}_x\otimes_\Z\prod_{\ell\neq p}\Z_\ell(1)
}$$
whose vertical arrows are the isomorphisms \eqref{eq_madonnadellolmo},
and whose top (resp. bottom) horizontal arrow is induced by
the natural morphism $X(\bar y)\to X$ (resp. by the specialization
map $\underline M{}_x\to\underline M{}_{\bar y}$).
Now, by corollary \ref{cor_normal-and-CM}, the scheme $U$
is connected and normal, hence the restriction functor
$\bCov(U)\to\bTame(X,(X,\underline M)_\tr)$ is fully faithful
(lemma \ref{lem_replace}(iii)); by \cite[Exp.V, Prop.6.9]{SGA1}
and theorem \ref{th_Abhyankar_log}, it follows that the induced
group homomorphism
\set\begin{equation}\label{eq_gimme-this}
\pi_1((X,\underline M)_\et,\xi)\to\pi_1(U_\et,\xi)
\end{equation}
is surjective (see \eqref{subsec_supply}). Clearly, the
image of $\phi_y$ lies in the kernel of \eqref{eq_gimme-this}.
Especially:
$$
(U,\underline M{}_{|U})_r\neq\emptyset\Rightarrow
\pi_1(U_\et,\xi)=\{1\}
$$
since, for any geometric point $\bar y$ of $(X,\underline M)_r$,
the specialization map $\underline M{}_x\to\underline M{}_{\bar y}$
is an isomorphism (notation of definition \ref{def_trivial-locus}(i)).
Thus, suppose that $(X,\underline M)_r=\{x\}$, set 
$P:=\underline M{}_x^\sharp$, and denote by $\psi:X\to\Spec\,P$
the natural continuous map; we have
$\underline M{}_{\bar y}^\sharp=P_{\psi(y)}$, and if we let
$F_y:=P\setminus\psi(y)$, we get a short exact sequence of
free abelian groups of finite rank :
\set\begin{equation}\label{eq_punctured-log}
0\to F^\gp_y\to\underline M{}_x^{\sharp\gp}\to
\underline M{}_{\bar y}^{\sharp\gp}\to 0.
\end{equation}
It is easily seen that $\psi$ is surjective; hence, for
any $\fp\in\Spec\,P$ of height $r-1$, pick $y\in\psi^{-1}(\fp)$.
With this choice, we have $F^\gp_y=\Z$; considering
$\eqref{eq_punctured-log}^\vee$, we deduce that $\pi_1(U_\et,\xi)$
is a quotient of $\hat\Z'(1):=\prod_{\ell\neq p}\Z_\ell(1)$.
More precisely, let $(\Spec\,P)_{r-1}$ be the set of prime
ideals of $P$ of height $r-1$; then we have :

\begin{theorem}\label{th_pi_1-punctured}
If\/ $(X,\underline M)_r=\{x\}$, we have a natural identification :
\set\begin{equation}\label{eq_goto-source}
\frac{P^\vee}{\sum_{\fp\in(\Spec\,P)_{r-1}}P_\fp^\vee}
\otimes_\Z\hat\Z'(1)\isom\pi_1(U_\et,\xi)
\end{equation}
and these two abelian groups are cyclic of finite order prime to $p$.
\end{theorem}
\begin{proof} The foregoing discussion already yields a natural
surjective map as stated; it remains only to check the injectivity,
and to verify that the source is a cyclic finite group.

However, since $d\geq 2$, and $(X,\underline M)_r=\{x\}$, corollary
\ref{cor_same-height}(i) implies that $r\geq 2$, in which case
$(\Spec\,P)_{r-1}$ must contain at least two distinct elements,
say $\fp$ and $\fq$. Let $F:=P\setminus\fp$; the image $H$ of
$(P_\fq)^{\sharp\gp\vee}$ in $F^{\gp\vee}$ is clearly a non-trivial
subgroup, hence $F^{\gp\vee}/H$ is a cyclic finite group, and then
the same holds for the source of \eqref{eq_goto-source}.

Next, it follows easily from lemma \ref{lem_replace}(iii)
that \eqref{eq_gimme-this} identifies $\pi_1(U_\et,\xi)$
with the quotient of $\pi_1((X,\underline M)_\et,\xi)$ by
the sum of the images of the maps $\phi_y$, for $y$ ranging
over all the points of $U$. However, it is clear that the
same sum is already spanned by the sum of the images of
the $\phi_y$ such that $y\in(X,\underline M)_{r-1}$; whence
the theorem.
\end{proof}

\begin{example} Let $P\subset\N^{\oplus 2}$ be the submonoid
of all pairs $(a,b)$ such that $a+b\in 2\N$. Clearly $P$ is
fine and saturated of dimension $2$. Let $K$ be any field;
the $K$-scheme $X:=\Spec\,K[P]$ is the singular quadric in
$\A^3_K$ cut by the equation $XY-Z^2=0$. It is easily seen
that $\Spec(K,P)_2$ consists of a single point $x$ (the vertex
of the cone); let $\bar x$ be a geometric point of $X$
localized at $x$, and $U\subset X(\bar x)$ the complement
of the closed point. Then $U$ is a normal $K$-scheme of dimension
$1$, and we have a natural isomorphism
$$
\pi_1(U,\xi)\isom\Z/2\Z.
$$
Indeed, a simple inspection shows that $P$ admits exactly two
prime ideals of height one, namely $\fp:=P\setminus(2\N\oplus\{0\})$
and $\fq:=P\setminus(\{0\}\oplus 2\N)$. Then,
$P_\fp=\{(a,b)\in\Z\oplus\N~|~a+b\in 2\N\}$, and similarly
$P_\fq$ is a submonoid of $\N\oplus\Z$. The quotients $P_\fp^\sharp$
and $P^\sharp_\fq$ are both isomorphic to $\N$, and are both generated
by the class of $(1,1)$. Let $\phi:P^\sharp_\fp\to\Z$ be a map
of monoids; then the image of $\phi$ in $P^{\gp\vee}$ is the
unique map of monoids $P\to\Z$ given by the rule : $(2,0)\mapsto 0$,
$(1,1)\mapsto\phi(1,1)$, $(0,2)\mapsto 2\phi(1,1)$. Likewise,
a map $\psi:P^\sharp_\fq\to\Z$ gets sent to the morphism $P\to\Z$
such that $(2,0)\mapsto 2\psi(1,1)$, $(1,1)\mapsto\psi(1,1)$, and
$(0,2)\mapsto 0$. We see therefore that
$(P_\fp)^{\sharp\gp\vee}+(P_\fq)^{\sharp\gp\vee}$ is a subgroup
of index two in $P^{\gp\vee}$, and the contention follows from
theorem \ref{th_pi_1-punctured}.
\end{example}

\subsection{Local acyclicity of smooth morphisms of log schemes}
\label{sec_acyclic-log}
In this section we consider a smooth and saturated morphism
$f:(X,\underline M)\to(Y,\underline N)$ of fine log schemes.
We fix a geometric point $\bar x$ of $X$, localized at a point
$x$, and let $\bar y:=f(\bar x)$. We shall suppose as well that
$Y$ is strictly local and normal, that $(Y,\underline N)_\tr$ is
a dense open subset of $Y$, and that $\bar y$ is localized
at the closed point $y$ of $Y$. Let $\bar\eta$ be a strict
geometric point of $Y$, localized at the generic point $\eta$;
to ease notation, set :
$$
U:=f^{-1}_{\bar x}(\eta)
\qquad
U_\tr:=(X(\bar x),\underline M(\bar x))_\tr\cap U
\qquad
\bar U:=U\times_{|\eta|}|\bar\eta|
\qquad
\bar U_\tr:=U_\tr\times_{|\eta|}|\bar\eta|
$$
and notice that $U_\tr$ is a dense open subset of $U$, by virtue of
proposition \ref{prop_generic-boundary}(ii,iv). Choose a geometric
point $\xi$ of $\bar U_\tr$, and let $\xi'$ be the image of $\xi$ in
$U_\tr$. There follows a short exact sequence of topological groups :
$$
1\to\pi_1(\bar U_{\tr,\et},\xi)\to\pi_1(U_{\tr,\et},\xi')
\xrightarrow{\ \alpha\ }\pi_1(|\eta|_\et,\bar\eta)\to 1.
$$
On the other hand, let $p$ be the characteristic exponent of
$\kappa:=\kappa(\bar y)$; for every integer $e$ such that $(e,p)=1$,
the discussion of \eqref{subsec_first-pairing} yields a
commutative diagram :
\set\begin{equation}\label{eq_many-eqs}
{\diagram
\underline N{}_{\bar y} \ar[r] \ar[d]_{\log f_{\bar x}} &
\kappa(\eta)^\times \ar[r] \ar[d] &
\Hom_\mathrm{cont}(\pi_1(|\eta|_\et,\bar\eta),\bmu_e(\kappa))
\ar[d]^{\Hom_\mathrm{cont}(\alpha,\bmu_e(\kappa))} \\
\underline M{}_{\bar x} \ar[r] & \Gamma(U_\tr,\cO^\times_U) \ar[r] &
\Hom_\mathrm{cont}(\pi_1(U_{\tr,\et},\xi'),\bmu_e(\kappa))
\enddiagram}
\end{equation}
such that the composition of the two top (rep. bottom) horizontal
arrows factors through $\underline N{}_{\bar y}^\sharp$ (resp.
$\underline M{}_{\bar x}^\sharp$). There follows a system of
natural group homomorphisms :
\set\begin{equation}\label{eq_thes-ame-dedu}
\Coker\,(\log f)^\gp_{\bar x}\to
\Hom_\mathrm{cont}(\pi_1(\bar U_{\tr,\et},\xi),\bmu_e(\kappa))
\qquad
\text{where $(e,p)=1$}
\end{equation}
which assemble into a group homomorphism :
\set\begin{equation}\label{eq_log-pairing}
\pi_1(\bar U_{\tr,\et},\xi)\to\Coker\,(\log f^\gp_{\bar x})^\vee
\otimes_\Z\prod_{\ell\neq p}\Z_\ell(1).
\end{equation}

\sset\subsubsection{}\label{subsec_relative-torsor}
We wish to give a second description of the map
\eqref{eq_log-pairing}, analogous to the discussion in
\eqref{subsec_Second-pairing}. To this aim, let $R$ be
a ring; for any monoid $P$, and any integer $e>0$, denote
by $\bee_P:P\to P$ the $e$-Frobenius map of $P$.
Let $\lambda:P\to Q$ be a local morphism of finitely generated
monoids. For any integer $e>0$, we get a commutative diagram
of monoids :
\set\begin{equation}\label{eq_pietiner}
{\diagram
P \ar[r]^-{\bee_P} \ar[d]_\lambda &
P \ar[d]^{\lambda_e} \ar[drr]^\lambda \\
Q \ar[r]^-{\mu'} & Q' \ar[rr]^{\bee_{Q|P}} & & Q
\enddiagram}
\end{equation}
whose square subdiagram is cocartesian, and such that
$\bee_{Q|P}\circ\mu'=\bee_Q$. The latter induces a
commutative diagram of log schemes :
$$
\xymatrix{
\Spec(R,Q) \ar[rr]^-{g_{Q|P}} \ar[drr]_g & &
\Spec(R,Q') \ar[r]^-{g'} \ar[d]_{g_e} &
\Spec(R,Q) \ar[d]^g \\
& & \Spec(R,P) \ar[r]^-{g_P} & \Spec(R,P)
}$$
whose square subdiagram is cartesian, and such that
$g_Q:=g'\circ g_{Q|P}$ is the morphism induced by $\bee_Q$.
Also, by example \ref{ex_fan-andlogscheme}, we have a commutative
diagram of monoidal spaces :
$$
\xymatrix{
\Spec(R,Q) \ar[rr]^-{g_{Q|P}} \ar[d]_{\psi_Q} & &
\Spec(R,Q') \ar[d]^{\psi_{Q'}} \ar[rr]^-{g_e} & & 
\Spec(R,P) \ar[d]^{\psi_P} \\
T_Q \ar[rr]^-\phi & & T_{Q'} \ar[rr] & & T_P
}$$
(where $\phi:=(\Spec\,\bee_{Q|P})^\sharp$) which determines
morphisms in the category $\cK$ :
$$
(\Spec(R,Q),T_Q,\psi_Q)\to(\Spec(R,Q'),T_{Q'},\psi_{Q'})\to
(\Spec(R,P),T_P,\psi_P).
$$

\sset\subsubsection{}\label{subsec_discuss-strict-loc}
Now, suppose that the closed point $\fm_Q$ of $T_Q$ lies in the
strict locus of $(\Spec\,\lambda)^\sharp$ (which just means that
$\lambda^\sharp$ is an isomorphism). Notice that the functor
$M\mapsto M^\sharp$ commutes with colimits (since it is a left
adjoint); taking into account lemma \ref{lem_localize-sharp},
we deduce that the square subdiagram of
$\eqref{eq_pietiner}^\sharp$ is still cocartesian, and therefore
$\bee^\sharp_{Q|P}$ is an isomorphism, so the same holds for $\phi$.

More generally, let $\fq\in\Spec\,Q$ be any prime ideal in the
strict locus of $(\Spec\,\lambda)^\sharp$; set $\fp:=\lambda^{-1}\fq$,
$\fr:=\phi(\fq)$, and recall that $T_{Q_\fq}:=\Spec\,Q_\fq$ is
naturally an open subset of $T_Q$ (see \eqref{subsec_height-in-T}).
A simple inspection shows that the restriction
$T_{Q_\fq}\to T_{Q'_\fr}$ of $\phi$ is naturally identified with
the morphism of affine fans $(\Spec\,\bee_{Q_\fq|P_\fp})^\sharp$.
Since $\fq$ is the closed point of $T_{Q_\fq}$, the foregoing
shows that $T_{Q_\fq}$ lies in the strict locus of $\phi$;
in other words, $\Str(\phi)$ is an open subset of $T_Q$, and
we have
\set\begin{equation}\label{eq_letargic}
\Str((\Spec\,\lambda)^\sharp)\subset\Str(\phi).
\end{equation}
Moreover, recall that $\Spec\,\bee_Q:T_Q\to T_Q$ is the
identity on the underlying topological space (see example
\ref{ex_multiply-by-n-in-fan}(i)), and by construction
it factors through $\phi$, so the latter is injective on
the underlying topological spaces. Especially,
$\Str(\phi)=\phi^{-1}\phi(\Str(\phi))$, and therefore
\set\begin{equation}\label{eq_strictloci}
\Str(g_{Q|P})=\psi_Q^{-1}(\Str(\phi))=
g^{-1}_{Q|P}(\psi^{-1}_{Q'}\phi(\Str(\phi))).
\end{equation}
Hence, set $g:=\Spec(R,\lambda)$; from \eqref{eq_letargic} and
\eqref{eq_strict-identity} we obtain :
$$
\Str(g)\subset\Str(g_{Q|P})
$$
and together with \eqref{eq_strictloci} we deduce that :
\begin{itemize}
\item
$\Str(g_{Q|P})$ is an open subset of $\Spec\,R[Q]$.
\item
$\psi^{-1}_{Q'}\phi(\Str(\phi))$ is a locally closed subscheme
of $\Spec\,R[Q']$.
\item
The restriction $\Str(g)\to\psi^{-1}_{Q'}\phi(\Str(\phi))$ of
$g_{Q|P}$ is a finite morphism.
\end{itemize}
Lastly, suppose that $e$ is invertible in $R$; in this case,
the morphisms $g_P$ and $g_Q$ are \'etale (proposition
\ref{prop_toric-smooth}), so the same holds for
$g_{Q|P}$ (corollary \ref{cor_sorite-smooth}(iii)).
From corollary \ref{cor_undercover}(i) we deduce that
the restriction
$\Str(g)\to\psi^{-1}_{Q'}\phi(\Str(\phi))$
of $g_{Q|P}$ is a (finite) \'etale covering.

\sset\subsubsection{}\label{subsec_Galois-relative}
In the situation of \eqref{subsec_relative-torsor}, suppose
additionally, that $R$ is a $\Z[1/e,\bmu_e]$-algebra (where
$\bmu_e$ is the $e$-torsion subgroup of $\C^\times$), $P$ and
$Q$ are fine monoids, and $\lambda$ is integral, so that $Q'$
is also fine. Set
$$
G_P:=\Hom_\Z(P^\gp,\bmu_e)
\qquad
G_Q:=\Hom_\Z(Q^\gp,\bmu_e)
\qquad
G_{Q|P}:=\Hom_\Z(\Coker\,\lambda^\gp,\bmu_e).
$$
Notice that the trivial locus $\Spec(R,P)_\tr$ is the open
subset $\Spec\,R[P^\gp]$, and likewise for $\Spec(R,Q)_\tr$
and $\Spec(R,Q')_\tr$; therefore, $\eqref{eq_pietiner}^\gp$
induces a cartesian diagram of schemes
$$
\xymatrix{
\Spec(R,Q')_\tr \ar[rr]^-{g'_\tr} \ar[d]_{g_{e,\tr}} & &
\Spec(R,Q)_\tr \ar[d]^{g_\tr} \\
\Spec(R,P)_\tr \ar[rr]^-{g_{P,\tr}} & & \Spec(R,P)_\tr.
}$$
Fix a geometric point $\tau'_Q$ of $\Spec(R,Q')_\tr$, and let
$\tau_Q:=g'(\tau'_Q)$, $\tau_P:=g(\tau_Q)$, $\tau'_P:=g_e(\tau'_Q)$.
It was shown in \eqref{subsec_Second-pairing} that $g^{-1}_P(\tau_P)$
is a $G_P$-torsor, so $g_{P,\tr}$ is a Galois \'etale covering,
corresponding to a continuous representation \eqref{eq_to-be-composed}
into $G_P$.
Hence $g^{\prime-1}(\tau_Q)$ is a $G_P$-torsor, and $g'_\tr$ is
a Galois \'etale covering, whose corresponding representation of
$\pi_1(\Spec(R,Q)_{\tr,\et},\tau_Q)$ is obtained by composing
\eqref{eq_to-be-composed} with the natural continuous
group homomorphism
$$
\pi_1(g_\tr,\tau_Q):\pi_1(\Spec(R,Q)_{\tr,\et},\tau_Q)\to
\pi_1(\Spec(R,P)_{\tr,\et},\tau_P).
$$
By the same token, $g_Q^{-1}(\tau_Q)$ is a $G_Q$-torsor, so
also $g_{Q,\tr}$ is a Galois \'etale covering, and a simple
inspection shows that the surjection
$$
g_Q^{-1}(\tau_Q)\to g^{\prime-1}(\tau_Q)
$$
induced by $g_{Q|P}$ is $G_Q$-equivariant, for the $G_Q$-action
on the target obtained from the map
$$
\Hom_\Z(\lambda^\gp,\bmu_e):G_Q\to G_P.
$$
The situation is summarized by the commutative diagram of continuous
group homomorphisms
$$
\xymatrix{
\pi_1(\Spec(R,Q')_{\tr,\et},\tau'_Q) \ar[r]
\ar[d]_{\pi_1(g_{e,\tr},\tau'_Q)} &
\pi_1(\Spec(R,Q)_{\tr,\et},\tau_Q) \ar[r]
\ar[d]_{\pi_1(g_\tr,\tau_Q)} &
G_Q \ar[d]^-{\Hom_\Z(\lambda^\gp,\bmu_e)} \\
\pi_1(\Spec(R,P)_{\tr,\et},\tau'_P) \ar[r] &
\pi_1(\Spec(R,P)_{\tr,\et},\tau_P) \ar[r] & G_P
}$$
whose horizontal right-most arrows are the maps \eqref{eq_to-be-composed}.
Consequently, $g^{-1}_{Q|P}(\tau'_Q)$ is a $G_{Q|P}$-torsor,
and $g_{Q|P,\tr}$ is a Galois \'etale covering, classified by
a continuous group homomorphism
\set\begin{equation}\label{classify}
\Ker\,\pi_1(g_{e,\tr},\tau'_Q)\to
\Ker\,\pi_1(g_\tr,\tau_Q)\to G_{Q|P}.
\end{equation}

\sset\subsubsection{}\label{subsec_strict-is-finite}
Let us return to the situation of \eqref{sec_acyclic-log}, and
assume additionally, that both $(X,\underline M)$ and
$(Y,\underline N)$ are fs log schemes. Take $R:=\cO_{Y,\bar y}$
in \eqref{subsec_relative-torsor}; by corollary
\ref{cor_charact-smoothness} and theorem
\ref{th_good-charts}(iii), we may assume that there exist 
\begin{itemize}
\item
a local and saturated morphism $\lambda:P\to Q$ of fine and
saturated monoids, such that $P$ is sharp, $Q^\times$ is
a free abelian group of finite type, say of rank $r$, and
$\Ass_\Z\,\Coker\,\lambda^\gp$ does not contain the characteristic
exponent of $\kappa(\bar y)$;
\item
a morphism of schemes $\pi:Y\to\Spec\,R[P]$, which is a section
of the projection $\Spec\,R[P]\to Y$, such that
$$
(Y,\underline N)=Y\times_{\Spec\,R[P]}\Spec(R,P)
\qquad
(X,\underline M)=Y\times_{\Spec\,R[P]}\Spec(R,Q)
$$
and $f$ is obtained by base change from the morphism
$g:=\Spec(R,\lambda)$. Moreover, the induced chart
$Q_X\to\underline M$ shall be local at the geometric point $\bar x$.
\end{itemize}
By claim \ref{cl_section-sharp}, we may further assume that the
projection $Q\to Q^\sharp$ admits a section $\sigma:Q^\sharp\to Q$,
such that $\lambda(P)$ lies in the image of $\sigma$. In this case,
$g$ factors through the morphism $\Spec(R,P)\to\Spec(R,Q^\sharp)$
induced by $\lambda^\sharp$, and $\sigma$ induces an isomorphism
of log schemes :
$$
\Spec(R,Q)=\G^{\oplus r}_{m,Y}\times_Y\Spec(R,Q^\sharp)
$$
(where $\G^{\oplus r}_{m,Y}$ denotes the standard torus of rank
$r$ over $Y$). In this situation, $X$ is smooth over
$Y\times_{\Spec\,R[P]}\Spec\,R[Q^\sharp]$, and more precisely 
\set\begin{equation}\label{eq_never-done}
(X,\underline M)=
\G^{\oplus r}_{m,Y}\times_{\Spec\,R[P]}\Spec(R,Q^\sharp).
\end{equation}
Summing up, after replacing $Q$ by $Q^\sharp$, we may assume
that $Q$ is also sharp, and \eqref{eq_never-done} holds with
$\Spec(R,Q)$ instead of $\Spec(R,Q^\sharp)$.
Moreover, the image of $x$ in $T_Q$ is the closed point $\fm_Q$.

Let $e>0$ be an integer which is invertible in $R$; by inspecting
the definition, it is easily seen that there exists a finite
separable extension $K_e$ of $\kappa(\eta)=\Frac(R)$, such that the
normalization $Y_e$ of $Y$ in $\Spec\,K_e$ fits into a commutative
diagram of schemes :
$$
\xymatrix{ Y_e \ar[r] \ar[d]_{\pi_e} & Y \ar[d]^\pi \\
           \Spec\,R[P] \ar[r]^-{g_P} & \Spec\,R[P]
}$$
(whose top horizontal arrow is the obvious morphism); namely,
$\pi$ is defined by some morphism of monoids $\beta:P\to R$,
and one takes for $K_e$ any subfield of $\kappa(\bar\eta)$
containing $\kappa(\eta)$ and the $e$-th roots of the elements
of $\beta(P)$. Set
$(Y_e,\underline N{}_e):=Y_e\times_{\Spec\,R[P]}\Spec(R,P)$,
and define log schemes $(X_e,\underline M{}_e)$,
$(X'_e,\underline M{}'_e)$ so that the two square subdiagrams
of the diagram of log schemes
$$
\xymatrix{
(X'_e,\underline M{}'_e) \ar[r]^{h_e} \ar[d] &
(X_e,\underline M{}_e) \ar[r] \ar[d] &
\G^{\oplus r}_{m,Y}\times_Y(Y_e,\underline N{}_e) \ar[d]^{\pi'_e} \\
\Spec(R,Q) \ar[r]^-{g_{Q|P}} &
\Spec(R,Q') \ar[r]^-{g_e} & \Spec(R,P)
}$$
are cartesian (here $\pi'_e$ is the composition of $\pi_e$
and the projection $\G^{\oplus r}_{m,Y}\times_YY_e\to Y_e$),
whence a commutative diagram :
\set\begin{equation}\label{eq_notations}
{\diagram (X'_e,\underline M{}'_e) \ar[r]^{h_e} \ar[rd]_{f'_e} &
(X_e,\underline M{}_e) \ar[d]^{f_e} \\
& (Y_e,\underline N{}_e).
\enddiagram}
\end{equation}
Notice that both $f_e$ and $f'_e$ are smooth and saturated morphisms
of fine log schemes. Also, by construction $Y_e$ is strictly local,
and $(Y_e,\underline N{}_e)_\tr$ is a dense subset of $Y_e$.
In other words, $f_e$ and $f'_e$ are still of the type considered
in \eqref{sec_acyclic-log}. Moreover, $h_e$ is \'etale, since
the same holds for $g_{Q|P}$, and the discussion in
\eqref{subsec_discuss-strict-loc} shows that the restriction
$$
\Str(f'_e)\to X_e
$$
is an \'etale morphism of schemes. Furthermore, the discussion
in \eqref{subsec_Galois-relative} shows that
$$
h_{e,\tr}:(X'_e,\underline M{}'_e)_\tr\to(X_e,\underline M{}_e)_\tr
$$
is a Galois \'etale covering.

\sset\subsubsection{}\label{subsec_no-duplicate}
More precisely, notice that $\lambda^\sharp=\log f_{\bar x}^\sharp$;
combining with \eqref{classify}, we deduce a continuous group
homomorphism :
\set\begin{equation}\label{eq_natur-continu}
\Ker\,\pi_1(f_{e,\tr},\xi'_e)\to
\Hom_\Z(\Coker\,(\log f)^\gp_{\bar x},\bmu_e(\kappa))
\end{equation}
where $\xi'_e$ is the image in $X_e$ of the geometric point $\xi$.
The geometric point $\bar y$ lifts uniquely to a geometric point
$\bar y_e$ of $Y_e$, localized at the closed point $y_e$, and the
pair $(\bar x,\bar y_e)$ determines a unique geometric point
$\bar x_e$ such that $f_e(\bar x_e)=\bar y_e$.
Also, since the field extension $\kappa(y)\to\kappa(y_e)$ is purely
inseparable, it is easily seen that the induced map
$f_e^{-1}(y_e)\to f^{-1}(y)$ is a homeomorphism; there follows
an isomorphism of $X(\bar x)$-schemes (\cite[Ch.IV, Prop.18.8.10]{EGA4}) :
$$
X_e(\bar x_e)\isom X(\bar x)\times_YY_e.
$$
Let $\eta_e$ be the generic point of $Y_e$; by construction,
$\bar\eta$ lifts to a geometric point $\bar\eta_e$ of $Y_e$,
localized at $\eta_e$, and we have continuous group homomorphisms :
\set\begin{equation}\label{eq_yahii}
\pi_1(\bar U_\tr,\xi)\isom
\Ker\,(\pi_1(U_{e,\tr},\xi'_e)\to\pi_1(\eta_e,|\bar\eta_e|))\to
\Ker\,\pi_1(f_{e,\tr},\xi'_e).
\end{equation}
The composition of \eqref{eq_natur-continu} and \eqref{eq_yahii}
is a continuous group homomorphism
\set\begin{equation}\label{eq_explict-quasi}
\pi_1(\bar U_\tr,\xi)\to
\Hom_\Z(\Coker\,(\log f)^\gp_{\bar x},\bmu_e(\kappa))
\end{equation}
whence, finally, a pairing
$$
\Coker\,(\log f)^\gp_{\bar x}\times\pi_1(\bar U_{\tr,\et},\xi)\to
\bmu_e(\kappa).
$$
We claim that this pairing agrees with the one deduced
from \eqref{eq_thes-ame-dedu}. Indeed, by tracing back
through the constructions, we see that \eqref{eq_explict-quasi}
is the homomorphism arising from the Galois covering of $\bar U_\tr$,
which is obtained from $g_{Q|P}$, after base change along the
composition
$$
\bar U_\tr\to X_e\to\Spec\,R[Q'].
$$
On the other hand, the discussion of \eqref{subsec_Second-pairing}
shows that the homomorphism $\pi_1(U_{\tr,\et},\xi')\to G_Q$
arising from the bottom row of \eqref{eq_many-eqs}, classifies
the Galois covering $C\to\bar U_\tr$ obtained from $g_Q$, by base
change along the same map. By the same token, the top row
of \eqref{eq_many-eqs} corresponds to the $G_P$-Galois covering
$C'\to|\eta|$ obtained by base change of $g_P$ along the
composition $|\eta|\to Y\to S_P$. The map $\log f_{\bar x}$
induces a morphism of schemes $C\to C'\times_{|\eta|}U_\tr$,
and \eqref{eq_thes-ame-dedu} corresponds to the $G_{Q|P}$-torsor
obtained from a fibre of this morphism. Evidently, this torsor
is isomorphic to $g_{Q|P}^{-1}(\tau'_Q)$, whence the contention.

\sset\subsubsection{}\label{subsec_Sigma}
In the situation of \eqref{sec_acyclic-log}, recall that
there is a natural bijection between the set of maximal
points of $f^{-1}_{\bar x}(\bar y)$, and the set $\Sigma$
of maximal points of the closed fibre of the induced map
\set\begin{equation}\label{eq_original-map}
\Spec\,\underline M{}_{\bar x}\to\Spec\,\underline N{}_{\bar y}
\end{equation}
(proposition \ref{prop_max-pts-fibre}). For every $\fq\in\Sigma$,
denote by $\eta_\fq$ the corresponding maximal point of
$f^{-1}_{\bar x}(\bar y)$, choose a geometric point
$\bar\eta_\fq$ localized at $\eta_\fq$, and let $X(\bar\eta_\fq)$
be the strict henselization of $X(\bar x)$ at $\bar\eta_\fq$.
Also, set
$$
U_\fq:=U\times_{X(\bar x)}X(\bar\eta_\fq)
\qquad
\bar U_\fq:=U_\fq\times_{|\eta|}|\bar\eta|
$$
and notice that $\bar U_\fq$ is an irreducible normal scheme.
Notice as well that $f$ induces a strict morphism
$(X(\bar\eta_\fq),\underline M(\bar\eta_\fq))\to(Y,\underline N)$
(theorem \ref{th_satura-smooth}(iii.a)), and therefore the
log structure of
$U_\fq\times_{X(\bar\eta_\fq)}(X(\bar\eta_\fq),\underline M(\bar\eta_\fq))$
is trivial.

Recall that $Z:=\bar U\setminus\bar U_\tr$ is a finite union
of irreducible closed subsets of codimension one in $\bar U$,
and $\cO_{\bar U,z}$ is a discrete valuation ring, for each 
maximal point $z\in Z$ (proposition \ref{prop_generic-boundary}(iv,v));
especially, the category $\bTame(\bar U,\bar U_\tr)$ is
well defined (definition \ref{def_tamely-ram}(iii)). We denote
$$
\bTame(f,\bar x)
$$
the full subcategory of $\bTame(\bar U,\bar U_\tr)$ consisting
of all the coverings $C\to\bar U_\tr$ such that, for every
$\fq\in\Sigma$, the induced covering
$$
C\times_{\bar U_\tr}\bar U_\fq\to\bar U_\fq
$$
is trivial ({\em i.e.} $C\times_{\bar U_\tr}\bar U_\fq$ is
a disjoint union of copies of $\bar U_\fq$).
It is easily seen that $\bTame(f,\bar x)$ is a Galois category
(see \cite[Exp.V, D\'ef.5.1]{SGA1}), and we obtain a fibre functor
for this category, by restriction of the usual fibre functor
$\phi\mapsto\phi^{-1}(\xi)$ defined on all \'etale coverings
$\phi$ of $\bar U_\tr$; we denote by $\pi_1(\bar U_\tr/Y_\et,\xi)$
the corresponding fundamental group. According to
\cite[Exp.V, Prop.6.9]{SGA1}, the fully faithful inclusion
$\bTame(f,\bar x)\to\bCov(\bar U_\tr)$ induces a continuous
surjective group homomorphism
\set\begin{equation}\label{eq_surge-pi}
\pi_1(\bar U_{\tr,\et},\xi)\to\pi_1(\bar U_\tr/Y_\et,\xi).
\end{equation}

\begin{proposition}\label{prop_surj-log-relat}
The map \eqref{eq_log-pairing} factors through \eqref{eq_surge-pi},
and the induced group homomorphism :
\set\begin{equation}\label{eq_compute-pi-tame}
\pi_1(\bar U_\tr/Y_\et,\xi)\to\Coker\,(\log f^\gp_{\bar x})^\vee
\otimes_\Z\prod_{\ell\neq p}\Z_\ell(1).
\end{equation}
is surjective.
\end{proposition}
\begin{proof} Let $\sigma_Y:Y\to Y^\mathrm{qfs}$ be the natural
morphism of schemes exhibited in remark \ref{rem_quasi-fine}(iv),
and set
$$
(Y,\underline N'):=
Y\times_{Y^\mathrm{qfs}}(Y,\underline N)^\mathrm{qfs}
\qquad
(X,\underline M'):=
(Y,\underline N')\times_{(Y,\underline N)}(X,\underline M).
$$
Since $f$ is saturated, both $(X,\underline M')$ and
$(Y,\underline N')$ are fs log schemes (see remark
\ref{rem_quasi-fine}(i)); also, the morphism of schemes
underlying the induced morphism of log schemes
$f':(X,\underline M')\to(Y,\underline N')$, agrees with
that underlying $f$. Moreover, by construction we have
$\underline N{}^{\prime\sharp}_{\bar z}=
(\underline N{}^\sharp_{\bar z})^\sat$
for every geometric point $\bar z$ of $Y$, especially
$(Y,\underline N')_\tr=(Y,\underline N)_\tr$. Likewise,
$\underline M{}^{\prime\sharp}_{\bar z}=
(\underline M{}^\sharp_{\bar z})^\sat$ (lemma
\ref{lem_little}(iii,iv)), therefore $\Str(f')=\Str(f)$, and
especially, $(X,\underline M')_\tr=(X,\underline M)_\tr$.
Furthermore, notice that the the natural map
\set\begin{equation}\label{eq_map-on-cokers}
\Coker\,(\log f)_{\bar x}^\gp\to\Coker\,(\log f')_{\bar x}^\gp
\end{equation}
is surjective, and its kernel is a quotient of
$(\underline M{}^\sat_{\bar x})^\times/\underline M{}^\times_{\bar x}$,
especially it is a torsion subgroup. However, the cokernel
of $(\log f)_{\bar x}^\gp$ equals the cokernel of
$(\log f^\sharp_{\bar x})^\gp$, hence it is torsion-free
(corollary \ref{cor_persist-integr}(ii)), so \eqref{eq_map-on-cokers}
is an isomorphism. Thus, we may replace $\underline N$ (resp.
$\underline M$) by $\underline N'$ (resp. $\underline M'$),
and assume from start that $f$ is a smooth, saturated morphism
of fs log schemes.

In this case, in light of the discussion of
\eqref{subsec_no-duplicate}, it suffices to prove that
\eqref{eq_explict-quasi} is a surjection, and that it factors
through $\pi_1(\bar U_\tr/Y_\et,\xi)$. To prove the surjectivity
comes down to showing that $\bar U_\tr\times_{X_e}X'_e$ is a
connected scheme.
However, let $x_e$ be the support of $\bar x_e$, and notice that
$\psi_P\circ\pi_e(y_e)=\fm_P$, the closed point of $T_P$. Since
$x$ maps to the closed point of $T_Q$, we deduce easily that the
image of $x_e$ in $T_{Q'}$ is the closed point $\fm_{Q'}$,
{\em i.e.} $x_e$ lies in the closed subscheme
$X_e\times_{S_Q'}\Spec\,\kappa\La Q'/\fm_{Q'}\Ra$. Notice
as well that $\bee_{Q|P}$ is of Kummer type (see definition
\ref{def_Kummer-monoids}); by proposition \ref{prop_crit-Kummer}(ii),
it follows that there exists a unique geometric point $\bar x{}'_e$
of $X'_e$ lying over $\bar x_e$, whence an isomorphism
of $X_e(\bar x_e)$-schemes (\cite[Ch.IV, Prop.18.8.10]{EGA4})
$$
X'_e(\bar x{}'_e)\isom X'_e\times_{X_e}X_e(\bar x_e).
$$
Hence $\bar U_\tr\times_{X_e}X'_e$ is an open subset of
$|\bar\eta_e|\times_{Y_e}X'_e(\bar x_e)$, and the latter
is an irreducible scheme, by proposition
\ref{prop_generic-boundary}(ii). We also deduce that the induced
morphism
$$
h_{e,\bar x}:(X'_e(\bar x{}'_e),\underline M{}'_e(\bar x{}'_e))
\to(X_e(\bar x_e),\underline M{}_e(\bar x_e))
$$
is a finite \'etale covering of log schemes. From the discussion
in \eqref{subsec_strict-is-finite}, we see that the restriction
of $h_{e,\bar x}$
$$
\Str(h_{e,\bar x})\to X_e(\bar x_e)
$$
is an \'etale morphism, and $\Str(h_{e,\bar x})$ contains the strict
locus of the induced morphism
$$
f'_{e,\bar x}:(X'_e(\bar x{}'_e),\underline M{}'_e(\bar x{}'_e))\to
(Y_e,\underline N{}_e).
$$
Since the field extension $\kappa(\bar y)\to\kappa(\bar y_e)$
is purely inseparable, $\bar\eta_\fq$ lifts uniquely to a geometric
point $\eta_{e,\fq}\in X_e(\bar x_e)$, and as usual we deduce
that the strict henselization $X_e(\bar\eta_\fq)$ of $X_e(\bar x_e)$
at $\bar\eta_{e,\fq}$ is isomorphic, as an $X_e(\bar x_e)$-scheme,
to $X_e(\bar x_e)\times_{X(\bar x)}X(\bar\eta_\fq)$. Moreover,
if $\eta_{e,\fq}$ is the support of $\bar\eta_{e,\fq}$, a simple
inspection shows that the fibre $h_{e,\bar x}^{-1}(\eta_{e,\fq})$
consists of maximal points of $f_{e,\bar x}^{\prime-1}(y_e)$.
By theorem \ref{th_satura-smooth}(iii.a), every point of
$h_{e,\bar x}^{-1}(\eta_{e,\fq})$ lies in $\Str(f'_{e,\bar x})$,
therefore the induced morphism
$X'_e\times_{X_e}X_e(\bar\eta_\fq)\to X_e(\bar\eta_\fq)$ is
finite and \'etale. Taking into account \eqref{subsec_supply},
we conclude that the \'etale covering
$X'_e\times_{X_e}\bar U_\tr\to\bar U_\tr$ is an object of
$\bTame(f,\bar x)$, whence the proposition.
\end{proof}

\sset\subsubsection{}\label{subsec_tameses}
Say that $Y=\Spec\,R$; for every algebraic field extension $K$ of
$\kappa(\eta)=\Frac(R)$, let $R_K$ be the normalization of $R$
in $K$, set $|\eta_K|:=\Spec\,K$ and
$$
Y_K:=\Spec\,R_K
\quad
(Y_K,\underline N{}_K):=Y_K\times_Y(Y,\underline N)
\quad
(X_K,\underline M{}_K):=Y_K\times_Y(X,\underline M).
$$
Moreover, let $f_K:(X_K,\underline M{}_K)\to(Y_K,\underline N{}_K)$
be the induced morphism, and $\bar y_K$ any geometric point
localized at the closed point $y_K$ of the strictly local scheme
$Y_K$; since the extension $\kappa(\bar y)\to\kappa(\bar y_K)$
is purely inseparable, there exists a unique geometric point
$\bar x_K$ of $X_K$ lifting $\bar x$, and we have an isomorphism
of $(X(\bar x),\underline M(\bar x))$-schemes :
$$
(X_K(\bar x_K),\underline M{}_K(\bar x{}_K))\isom
X_K\times_X(X(\bar x),\underline M(\bar x))=
Y_K\times_Y(X(\bar x),\underline M(\bar x)).
$$
Clearly the morphism $f_K$ is again of the type considered
in \eqref{subsec_Sigma}; especially, the maximal points of
$f^{-1}_{K,\bar x_K}(\bar y_K)$ are in natural bijection
with the elements of $\Sigma$, and it is natural to denote
$$
U_K:=U\times_Y|\eta_K|
\qquad
U_{K,\tr}:=U_\tr\times_Y|\eta_K|
\qquad
U_{K,\fq}:=U_{K,\tr}\times_{X(\bar x)} X(\bar\eta_\fq)
$$
for every $\fq\in\Sigma$. Then, let
$\bar\eta_{K,\fq}$ be the unique geometric point of
$f^{-1}_{K,\bar x_K}(\bar y_K)$ lying over $\bar\eta_\fq$,
and $\eta_{K,\fq}$ the support of $\bar\eta_{K,\fq}$; as usual,
we have
\set\begin{equation}\label{eq_from_hell}
X_K(\bar\eta_{K,\fq})=X(\bar\eta_\fq)\times_YY_K
\end{equation}
hence the above notation is consistent with the one introduced
for the original morphism $f$. Furthermore, if $z$ is any
maximal point of $\bar U\!\setminus\!\bar U_\tr$, the image
$z_K$ of $z$ in $U_K$ is a maximal point of
$U_K\!\setminus\!U_{K,\tr}$, and since the induced map
$\cO_{\!U_K,z_K}\to\cO_{\!\bar U,z}$ is faithfully flat,
proposition \ref{prop_generic-boundary}(iv,v) easily implies
that $\cO_{\!U_K,z_K}$ is a discrete valuation ring. We may
then denote
$$
\bTame(f,\bar x,K)
$$
the full subcategory of $\bTame(U_K,U_{K,\tr})$, consisting
of those objects  $C\to U_{K,\tr}$, such that the induced
covering $C\times_{U_{K,\tr}}U_{K,\fq}\to U_{K,\fq}$ is
trivial, for every $\fq\in\Sigma$. We have a natural functor
\set\begin{equation}\label{eq_tameses}
\Pscolim{K}\bTame(f,\bar x,K)\to\bTame(f,\bar x)
\end{equation}
where the $2$-colimit ranges over the filtered family
of all finite separable extensions $K$ of $\kappa(\eta)$.

\begin{lemma}\label{lem_descend-tameses}
The functor \eqref{eq_tameses} is an equivalence.
\end{lemma}
\begin{proof} Let $\bar h:\bar C\to\bar U_\tr$ be an
object of the category $\bTame(f,\bar x)$. According to
\cite[Ch.IV, Prop.17.7.8(ii)]{EGA4} and
\cite[Ch.IV, Prop.8.10.5]{EGAIV-3}, we may find a finite
extension $K$ of $\kappa(\eta)$, such that $\bar h$ descends
to a finite \'etale morphism
$$
h_K:C_K\to U_{K,\tr}.
$$
Let $C'_K$ (resp. $\bar C{}'$) denote the normalization of
$C_K$ (resp. of $\bar C$) over $U_K$ (resp. over $\bar U$).
Since the morphism $|\bar\eta|\to|\eta_K|$ is ind-\'etale,
we have $\bar C{}'=C'_K\times_{|\eta_K|}|\bar\eta|$
(\cite[Ch.IV, Prop.17.5.7]{EGA4}), and it follows easily that
$C_K$ is tamely ramified along the divisor
$U_K\!\setminus\!U_{K,\tr}$.
From \eqref{eq_from_hell} we get a natural isomorphism :
$$
\bar U_\fq\isom U_{K,\fq}\times_{|\eta_K|}|\bar\eta|.
$$
Thus, after replacing $K$ by a larger finite separable extension
of $\kappa(\eta)$, we may assume that the induced morphism
$C_K\times_{U_{K,\tr}}U_{K,\fq}\to U_{K,\fq}$ is a trivial
\'etale covering, for every $\fq\in\Sigma$.

This shows that \eqref{eq_tameses} is essentially surjective;
likewise one shows the full faithfulness : the details shall
be left to the reader.
\end{proof}

\sset\subsubsection{}\label{subsec_F_queue-II}
In the situation of \eqref{subsec_tameses}, let $K$ be an
algebraic extension of $\kappa(\eta)$, and
$$
h:C\to U_{K,\tr}
$$
any object of $\bTame(U_K,U_{K,\tr})$, and denote by $C'$ the
normalization of $X_K(\bar x_K)$ in $C$. We claim that there
exists a largest non-empty open subset
$$
E(h)\subset X_K(\bar x_K)
$$
such that the restriction $h^{-1}E(h)\to E(h)$ of $h$ is \'etale.
Indeed, in any case, $h$ restricts to an \'etale morphism on a
dense open subset containing $U_{K,\tr}$, and there exists
a largest open subset $E'\subset C'$ such that $h_{|E'}$
is \'etale (claim \ref{cl_trivial-Zorni} and lemma
\ref{lem_replace}(iii)). Then it is easily seen that
$E(h):=X_K(\bar x_K)\!\setminus\!h(C'\!\setminus\!E')$ will do.

\begin{lemma}\label{lem_F_queue}
With the notation of \eqref{subsec_F_queue-II}, the category
$\bTame(f,\bar x,K)$ is the full subcategory of
$\bTame(U_K,U_{K,\tr})$ consisting of those objects
$h:C\to U_{K,\tr}$ such that $E(h)$ contains the maximal
points of $X_K(\bar x_K)\times_{Y_K}|y_K|$.
\end{lemma}
\begin{proof} In view of claim \ref{cl_include-pts}, this
characterization is a rephrasing of the definition of the
category $\bTame(f,\bar x,K)$.
\end{proof}

\sset\subsubsection{}\label{subsec_F_queue}
Keep the notation of \eqref{subsec_F_queue-II}, and suppose
that $h$ is an object of $\bTame(f,\bar x,K)$. Fix $\fq\in\Sigma$;
then lemma \ref{lem_F_queue} says that $\eta_{K,\fq}\in E(h)$.
Thus, we obtain a functor
$$
\bTame(f,\bar x,K)\to\bCov(|\eta_{K,\fq}|)
\qquad :\qquad
C\mapsto C'\times_{X_K(\bar x_K)}|\eta_{K,\fq}|.
$$
However, the natural morphism $|\eta_{K,\fq}|\to|\eta_\fq|$
is radicial, hence it induces an equivalence
$$
\bCov(|\eta_\fq|)\isom\bCov(|\eta_{K,\fq}|)
$$
(lemma \ref{lem_replace}(i)). Combining these two functors
in the special special case where $K:=\kappa(\bar\eta)$,
we get  a functor
\set\begin{equation}\label{eq_special_tameses}
\bTame(f,\bar x)\to\bCov(|\eta_\fq|)
\quad :\quad
(C\to\bar U_\tr)\mapsto(C_{|\eta_\fq}\to|\eta_\fq|).
\end{equation}
Now, the rule $\phi\mapsto\phi^{-1}(\bar\eta_\fq)$ yields a fibre
functor for the Galois category $\bCov(|\eta_\fq|)$;
by composition with \eqref{eq_special_tameses}, we deduce a fibre
functor for $\bTame(f,\bar x)$, whose group of automorphisms we
denote $\pi_1(\bar U_\tr/Y_\et,\bar\eta_\fq)$. Also, set
$F(\fq):=\underline M{}_{\bar x}\!\setminus\!\fq$, and notice
that the structure map
$\underline M{}_{\bar x}\to\cO_{X(\bar x),\bar x}$ induces a
group homomorphism
\set\begin{equation}\label{eq_should-I}
F(\fq)^\gp\to\kappa(\eta_\fq)^\times.
\end{equation}

\begin{lemma}\label{lem_horse-with-no-name}
With the notation of \eqref{subsec_F_queue}, we have :
\begin{enumerate}
\item
The natural map $\Coker\,(\log f_{\bar x})\to F(\fq)$
induces a commutative diagram of groups
$$
\xymatrix{
\pi_1(|\eta_\fq|_\et,\bar\eta_\fq) \ar[r] \ar[d]_\alpha
& \pi_1(\bar U_\tr/Y_\et,\bar\eta_\fq) \ar[r]
& \pi_1(\bar U_\tr/Y_\et,\xi) \ar[d]^\beta \\
F(\fq)^{\gp\vee}\otimes_\Z\prod_{\ell\neq p}\Z_\ell(1)
\ar[rr]^-\gamma & & \Coker\,(\log f^\gp_{\bar x})^\vee
\otimes_\Z\prod_{\ell\neq p}\Z_\ell(1)
}$$
where $\beta$ is \eqref{eq_compute-pi-tame}, and $\alpha$
is deduced from \eqref{eq_should-I}, as in the discussion
of \eqref{subsec_first-pairing}.
\item
$\alpha$ is surjective, and $\gamma$ is an isomorphism.
\end{enumerate}
\end{lemma}
\begin{proof}(i): The proof amounts to unwinding the definitions,
and shall be left as an exercice for the reader. Notice that
the second arrow on the top row is only well-defined up
to inner automorphisms, but since the groups on the bottom
row are abelian, the ambiguity does not affect the statement.

(ii): Notice that $\log f_{\bar x}$ restricts to a map of
monoids $\cO^\times_{Y,\bar y}\to F(\fq)$, which induces
an isomorphism $\Coker(\log f)^\sharp\isom F(\fq)^\sharp$;
we deduce that $\gamma$ is an isomorphism. Next, let $Z$
be the topological closure of $\eta_\fq$ in $X(\bar x)$,
and endow $Z$ with its reduced subscheme structure; set
also $(Z,\underline M(Z)):=Z\times_{X(\bar x)}(X,\underline M)$.
The map $\alpha$ factors as a composition
$$
\pi_1(|\eta_\fq|_\et,\bar\eta_\fq)\to
\pi_1((Z,\underline M(Z))_{\tr,\et},\bar\eta_\fq)\to
F(\fq)^{\gp\vee}\otimes_\Z\prod_{\ell\neq p}\Z_\ell(1)
$$
where the first map is surjective, by lemma
\ref{lem_replace}(ii). Lastly, notice that
$\underline M(Z)_{\red,\bar x}=F(\fq)_\circ$; by propositions
\ref{prop_fund-cuneo} and \ref{prop_max-pts-fibre}(ii), it
follows that the second map is surjective as well, so the
proof of (ii) is complete.
\end{proof}

\sset\subsubsection{}\label{subsec_tame-to-inf}
In the situation of \eqref{subsec_annoying-gen-II}, suppose that
$Y_i$ is a strictly local normal scheme for every $i\in I$, and
the transition morphisms $Y_j\to Y_i$ are local and dominant,
for every morphism $i\to j$ in $I$. Let $\bar x$ be a geometric
point of $X$, and denote by $\bar x_i$ the image of $\bar x$ in
$X_i$, for every $i\in I$. Suppose that the image $\bar y$ of
$\bar x$ in $Y$ is localized at the closed point. Also, let $\bar\eta$
be a strict geometric point of $Y$, localized at the generic
point $\eta_i$, and denote by $\bar\eta_i$ (resp. $\bar y_i$)
the strict image of $\bar\eta$ (resp. $\bar y$) in $Y_i$ (see
definition \ref{def_strict-loc}(v)).

\begin{lemma}\label{lem_descend-smooth-sat}
In the situation of \eqref{subsec_tame-to-inf}, suppose that
$(g,\log g):(X,\underline M)\to(Y,\underline N)$ is a smooth
and saturated morphism of fine log schemes.
Then there exist $i\in I$, and a smooth and saturated morphism
$(g_i,\log g_i):(X_i,\underline M{}_i)\to(Y_i,\underline N{}_i)$
of fine log schemes, such that $\log g=\pi_i^*\log g_i$.
\end{lemma}
\begin{proof} By corollary \ref{cor_smooth-charct}, we can
descend $(g,\log g)$ to a smooth morphism $(g_i,\log g_i)$
of fine log schemes, and after replacing $I$ by $I/i$, we
may assume that $i=0$. Then the contention follows from
corollary \ref{cor_descend-chart-from-infty}(ii).
\end{proof}

\sset\subsubsection{}
Keep the situation of \eqref{subsec_tame-to-inf}, and suppose that
$(g_0,\log g_0):(X_0,\underline M{}_0)\to(Y_0,\underline N{}_0)$
is a smooth and saturated morphism of fine log schemes; set
$$
(X_i,\underline M{}_i):=X_i\times_{X_0}(X_0,\underline M{}_0)
\qquad
(Y_i,\underline N{}_i):=Y_i\times_{Y_0}(Y_0,\underline Y{}_0)
$$
and denote
$(g_i,\log g_i):(X_i,\underline M{}_i)\to(Y_i,\underline N{}_i)$
the induced morphism of log schemes, for every $i\in I$.
Also, let $(g,\log g):(X,\underline M)\to(Y,\underline N)$
be the limit of the system of morphisms $((g_i,\log g_i)~|~i\in I)$.
These are morphisms of the type considered in \eqref{sec_acyclic-log},
so we may define $U_i:=g^{-1}_{i,\bar x_i}(\eta_i)$, and introduce
likewise the schemes $U_{i,\tr}$, $\bar U_i$ and $\bar U_{i,\tr}$
as in \eqref{sec_acyclic-log}. Moreover, set
$Z_i:=\bar U_i\!\setminus\!\bar U_{i,\tr}$ for every $i\in I$;
clearly $Z_i=Z_j\times_{X_j(\bar x_j)}X_i(\bar x_i)$ for every
morphism $i\to j$ in $I$. Also, each $Z_i$ is a finite union
of irreducible subsets of codimension one, and for every $i\to j$
in $I$, the transition morphisms $X_i\to X_j$ restrict to maps
$$
\Max\,Z_i\to\Max\,Z_j
\qquad
\Max\,X_i(\bar x_i)\times_{Y_i}|\bar y_i|\to
\Max\,X_j(\bar x_j)\times_{Y_j}|\bar y_j|
$$
Combining proposition \ref{prop_tame-infty}(ii) and lemma
\ref{lem_F_queue}, we deduce a fully faithful functor
\set\begin{equation}\label{eq_almost-donut}
\Pscolim{i\in I}\bTame(g_i,\bar x_i)\to\bTame(g,\bar x).
\end{equation}

\begin{lemma}\label{lem_chat_invisible}
The functor \eqref{eq_almost-donut} is an equivalence.
\end{lemma}
\begin{proof} It remains only to show the essential surjectivity.
Hence, let $h$ be a given object of $\bTame(g,\bar x)$; by
proposition \ref{prop_tame-infty}(ii), we know that there exists
$j\in I$ such that $h$ descends to an \'etale  covering
$h_j:V_j\to\bar U_{j,\tr}$, tamely ramified along $Z_j$, and after
replacing $I$ by $I/i$, we may assume that $j$ is the final object
of $I$, and define $h_i:=\bar U_{i,\tr}\times_{\bar U_{j,\tr}}h_j$
for every $i\to j$ in $I$. Now, let $E'\subset E(h)$ be a
constructible open subset containing the maximal points of
$X(\bar x)\times_Y|\bar y|$. For every $i\in I$, let $\bar Y_i$
be the normalization of $Y_i$ in $\Spec\,\kappa(\bar\eta_i)$,
and set $\bar X_i:=X_i(\bar x_i)\times_{Y_i}\bar Y_i$; according
to \cite[Ch.IV, Th.8.3.11]{EGAIV-3}, there exists $i\in I$ such
that $E'$ descends to a constructible open subset $E'_i\subset\bar X_i$,
and then necessarily $E'_i$ contains all the maximal points of
$X_i(\bar x_i)\times_{Y_i}|\bar y_i|$. As usual, we may assume
that $i$ is the final object, so $E'_i$ is defined for every
$i\in I$. Lastly, since $h$ extends to an \'etale covering on
$E'$, we see that $h_i$ extends to an \'etale covering of $E'_k$,
for some $k\in I$ (\cite[Ch.IV, Prop.17.7.8(i)]{EGA4}). In view
of lemma \ref{lem_F_queue}, the contention follows.
\end{proof}

\begin{theorem}\label{th_acyclicity-log-smooth}
The map \eqref{eq_compute-pi-tame} is an isomorphism.
\end{theorem}
\begin{proof} Arguing as in the proof of proposition
\ref{prop_surj-log-relat}, we may assume that both
$(X,\underline M)$ and $(Y,\underline N)$ are fs log schemes,
and in view of proposition \ref{prop_surj-log-relat},
we need only show that \eqref{eq_compute-pi-tame} is injective.
This comes down to the following assertion. For every object
$\bar h:\bar C\to\bar U_\tr$ of the category $\bTame(f,\bar x)$,
the induced action of $\pi_1(\bar U_\tr,\xi)$ on $\bar h{}^{-1}(\xi)$
factors through the quotient
$\Coker\,(\log f^\gp_{\bar x})^\vee\otimes_\Z\prod_{\ell\neq p}\Z_\ell(1)$.

$\bullet$\ \
By lemma \ref{lem_descend-tameses}, there exist a finite separable
extension $K$ of $\kappa(\eta)$, and an object $h:C_K\to U_{K,\tr}$
of $\bTame(f,\bar x,K)$, with an isomorphism
$C_K\times_{U_{K,\tr}}\bar U_\tr\isom\bar C$ of $\bar U_\tr$-schemes.
Since $\log f_{\bar x}=\log f_{K,\bar x_K}$, the theorem
will hold for the morphism $f$ and the point $\bar x$, if
and only if it holds for $f_K$ and the point $\bar x_K$.

\begin{claim}\label{cl_dim-one-ok}
The theorem holds if $Y$ is noetherian of dimension one.
\end{claim}
\begin{pfclaim} In this case, $Y$ is the spectrum of a strictly
henselian discrete valuation ring $R$, and the same then holds
for $Y_K$. Hence, we may replace throughout $f$ by $f_K$, and
assume from start that there exists an object $h:C\to U_\tr$
of $\bTame(f,\bar x,\kappa(\eta))$, with an isomorphism
$\bar C\isom C\times_{U_\tr}\bar U_\tr$ of $\bar U_\tr$-schemes. 
Endow $Y$ with the fine log structure $\underline N'$ such that
$\Gamma(Y,\underline N')=R\setminus\!\{0\}$; since
$(Y,\underline N)_\tr$ is dense in $Y$, we have a well defined
morphism of log schemes
$\pi:(Y,\underline N')\to(Y,\underline N)$, which is the
identity on the underlying schemes. Set $(X,\underline M'):=
(Y,\underline N')\times_{(Y,\underline N)}(X,\underline M)$.
Then $(Y,\underline N')$ is a regular log scheme, and consequently
the same holds for $(X,\underline M')$, by theorem
\ref{th_smooth-preserve-reg}.

Furthermore, $\pi$ trivially restricts to a strict morphism on
the open subset $|\eta|$, hence the induced morphism
$(X,\underline M')\times_Y|\eta|\to(X,\underline M)\times_Y|\eta|$
is an isomorphism, especially $U_{\tr}$ is the trivial locus
of $(X(\bar x),\underline M'(\bar x))\times_Y|\eta|$.
However, it is easily seen that
$(X(\bar x),\underline M'(\bar x))_\tr$ does not
intersect the closed fibre $f_{\bar x}^{-1}(\bar y)$,
so finally $U_\tr=(X(\bar x),\underline M'(\bar x))_\tr$.

From theorem \ref{th_Abhyankar_log}, we deduce that $h$ extends
to an \'etale covering of $(X(\bar x),\underline M'(\bar x))$.
Then, arguing as in \eqref{sec_acyclic-log} we get a
commutative diagram of groups :
$$
\xymatrix{ \pi_1(\bar U_{\tr,\et},\xi) \ar[r] \ar[d] &
\Coker\,(\log f^\gp_{\bar x})^\vee
\otimes_\Z\prod_{\ell\neq p}\Z_\ell(1) \ar[d] \\
\pi_1(U_{\tr,\et},\xi') \ar[r]^-\alpha &
\underline M{}^{\gp\vee}_{\bar x}\otimes_\Z\prod_{\ell\neq p}\Z_\ell(1)
}$$
whose top horizontal arrow is \eqref{eq_log-pairing}, and
whose right vertical arrow is deduced from the projection
$\underline M{}^{\gp}_{\bar x}\to\Coker\,(\log f^\gp_{\bar x})$.
Lastly, proposition \ref{prop_fund-cuneo} shows that the natural
action of $\pi_1(U_{\tr,\et},\xi')$ on $h^{-1}(\xi')$
factors through $\alpha$, so the claim follows.
\end{pfclaim}

$\bullet$\ \
Next, suppose that $Y$ is an arbitrary normal, strictly
local scheme. The discussion in \eqref{subsec_no-duplicate}
implies that, in order to prove the theorem, it suffices to
find an integer $e>0$, a $(X,\underline M)$-scheme
$(X'_e,\underline M{}'_e)$ as in \eqref{subsec_strict-is-finite},
and a geometric point $\bar x{}'_e$ of $X'_e$ lying over $\bar x$,
such that $h\times_{X(\bar x)}X'_e(\bar x{}'_e)$ is a trival
covering.  To this aim, we write $Y$ as the limit of a cofiltered
system $(Y_i~|~i\in I)$ of strictly local excellent and
normal schemes (lemma \ref{lem_excell}), and we denote by
$\eta_i$ the generic point of $Y_i$, for every $i\in I$. By lemma
\ref{lem_descend-smooth-sat}, we may then descend $(f,\log f)$
to a smooth and saturated morphism
$(f_i,\log f_i):(X_i,\underline M{}_i)\to(Y_i,\underline N{}_i)$,
for some $i\in I$, and as usual, we may assume that $i$ is
the final object of $I$. Let $\bar x_i$ be the image of
$\bar x$ in $X_i$; by lemmata \ref{lem_chat_invisible} and
\ref{lem_descend-tameses}, the object $h$ of
$\bTame(f,\bar x,\kappa(\eta))$ descends to an object $h_i$
of $\bTame(f_i,\bar x_i,K)$, for some $i\in I$, and some finite
separable extension $K$ of $\kappa(\eta_i)$. It suffices therefore
to find $e>0$, a
$(X_i,\underline M{}_i)$-scheme $(X'_{i,e},\underline M{}'_{i,e})$,
and a geometric point $\bar x{}'_{i,e}$ of $X'_{i,e}$ lying
over $\bar x_i$, such that
$h_i\times_{X_i(\bar x_i)}X'_{i,e}(\bar x{}'_{i,e})$
is a trivial covering. In other words, we may replace throughout
$Y$ by $Y_{i,K}$, and assume from start that $Y$ is excellent,
and $\bar h$ descends to an object $h:C\to U_\tr$ of
$\bTame(f,\bar x,\kappa(\eta))$.

$\bullet$\ \
By \cite[Ch.II, Prop.7.1.7]{EGAII} (see also remark
\ref{rem_correction-of-EGAII}) we may find a discrete valuation
ring $\sV$ and a local injective morphism $R\to\sV$ inducing an
isomorphism on the respective fields of fractions.
Let $\sV^\sh$ be the strict henselization of $\sV$
(at a geometric point whose support is the closed point),
and set
$$
\sY:=\Spec\,\sV^\sh
\quad
(\sY,\underline\sN):=\sY\times_Y(Y,\underline N)
\quad
(\sX,\underline\sM):=\sY\times_Y(X,\underline M).
$$
Also, let $\sff:(\sX,\underline\sM)\to(\sY,\underline\sN)$
be the induced morphism.
Denote by $\bar\sy$ a geometric point localized at the closed
point $\sy$ of $\sY$; also, pick any geometric point $\bar\sx$ of
$\sX$, whose image in $X$ is $\bar x$; the induced morphism
\set\begin{equation}\label{eq_base-change-V}
(\sX(\bar\sx),\underline\sM(\bar\sx))\to
(X(\bar x),\underline M(\bar x))
\end{equation}
restricts to a flat morphism
$\sff^{-1}_{\bar\sx}(\bar\sy)\to f^{-1}_{\bar x}(\bar y)$
and from proposition \ref{prop_max-pts-fibre}, we see that
the latter induces a bijection between the sets of maximal
points of the two fibres. On the other hand, let $\bar\eta_\sV$
denote a geometric point localized at the generic point $\eta_\sV$
of $\sY$; then \eqref{eq_base-change-V} restricts to an ind-\'etale
morphism $\sff^{-1}_{\bar\sx}(\eta_\sV)\to f^{-1}_{\bar x}(\eta)$.
Hence, set
$$
\sU_\tr:=
(\sX(\bar\sx),\underline\sM(\bar\sx))_\tr\times_{\sY}|\eta_\sV|.
$$
From lemma \ref{lem_F_queue}, it follows easily that the
covering $C\times_{U_\tr}\sU_\tr\to\sU_\tr$ is an object
of $\bTame(\sff,\bar\sx,\kappa(\eta_\sV))$.

For any integer $e>0$ invertible in $R$, pick a
$(Y,\underline N)$-scheme $(Y_e,\underline N{}_e)$ as in
\eqref{subsec_strict-is-finite}, so that we may
define the \'etale morphism
$(X'_e,\underline M{}'_e)\to(X_e,\underline M{}_e)$
of $(Y_e,\underline N{}_e)$-schemes as in \eqref{eq_notations}.
Notice that the morphism
$(X'_e,\underline M{}'_e)\to(Y_e,\underline N{}_e)$ is again of
the type considered in \eqref{sec_acyclic-log}, and there
exists, up to isomorphism, a unique geometric point
$\bar x{}'_e$ of $X'_e$ lifting $\bar x$; moreover, for any
geometric point $\bar y_e$ supported at $y_e$, the induced
map
$$
\Spec\,\underline M{}'_{e,\bar x{}'_e}\to
\Spec\,\underline N{}'_{e,\bar y{}_e}
$$
is naturally identified with \eqref{eq_original-map}.
Likewise, pick a $(\sY,\underline\sN)$-scheme
$(\sY_e,\underline\sN{}_e)$ in the same fashion, and denote
by $\eta_e$ (resp. $\eta_{\sV,e}$) the generic point of $Y_e$
(resp. of $\sY_e$), and by $\sy_e\in\sY_e$ (resp. $y_e\in Y_e$)
the closed point. We may choose $\sY_e$ so that
$\kappa(\eta_{\sV,e})$ contains $\kappa(\eta_e)$, in which
case we have a strict morphism
$$
(\sY_e,\underline\sN{}_e)\to(Y_e,\underline N{}_e)
$$
of log schemes, and we may set $(\sX'_e,\underline\sM{}'_e):=
\sY_e\times_{Y_e}(X'_e,\underline M{}'_e)$. Again, there
exists, up to isomorphism, a unique geometric point
$\bar\sx{}'_e$ of $\sX'_e$ lifting $\bar\sx$, and by claim
\ref{cl_dim-one-ok}, we may assume that both $e$ and
$\kappa(\eta_{\sV,e})$ have been chosen large enough, so that
the base change
$$
h\times_{X(\bar x)}\sX'_e(\bar\sx{}'_e)
$$
shall be a trivial \'etale covering. Hence, we may replace
$Y$ by $Y_e$, $\sY$ by $\sY_e$, $X$ by $X'_e$, and $h$ by
$h\times_{X(\bar x)}X'_e(\bar x{}'_e)$, and assume
from start that $C\times_{U_\tr}\sU_\tr$ is a trivial covering
of $\sU_\tr$. The theorem will follow, once we show that -- in
this case -- $h$ is a trivial \'etale covering.

Let $C'$ (resp. $\sC'$) be the normalization of $X(\bar x)$
in $C$ (resp. of $\sX(\bar\sx)$ in $C\times_{U_\tr}\sU_\tr$),
and define $E(h)$ as in \eqref{subsec_F_queue-II}.
Then $\sC'$ is a trivial \'etale covering of $\sX(\bar\sx)$,
and since $X(\bar x)$ is excellent, the induced morphism
$h':C'\to X(\bar x)$ is finite.

\begin{claim}\label{cl_new-Eprime}
For every maximal point $\eta_\fq$ of
$f^{-1}_{\bar x}(\bar y)$, the induced covering
$C_{|\eta_\fq}\to|\eta_\fq|$ is trivial (notation of
\eqref{eq_special_tameses}).
\end{claim}
\begin{pfclaim} Let $Z_\fq$ denote the topological closure
of $\{\eta_\fq\}$ in $X(\bar x)$, and endow $Z_\fq$ with
its reduced subscheme structure; then $E_\fq:=E(h)\cap Z_\fq$
is non-empty (lemma \ref{lem_F_queue}), and geometrically
normal (proposition \ref{prop_max-pts-fibre}(ii) and corollary
\ref{cor_normal-and-CM}). Also, \eqref{eq_base-change-V} induces
an isomorphism of $\kappa(\sy)$-schemes
$$
E_\fq\times_{X(\bar x)}\sX(\bar x)\isom
E_\fq\times_{\kappa(y)}\kappa(\sy).
$$
Moreover, $Z_\fq$ is strictly local, and
$Z_\fq\times_{X(\bar x)}\sX(\bar\sx)$ is the strict
henselization of $Z_\fq\times_{\kappa(y)}\kappa(\sy)$
at the point $\sx$. Furthermore, the morphism
$h'':=h'\times_{X(\bar x)}E_\fq$ is an \'etale covering
of $E_\fq$, and $h''\times_{\kappa(y)}\kappa(\sy)$
is naturally identified with the restriction of $\sC'$
to the subscheme $E_\fq\times_{X(\bar x)}\sX(\bar\sx)$
(\cite[Ch.IV, Prop.17.5.7]{EGA4}), hence it is a trivial
covering. By example \ref{ex_bcov}, it follows that $h''$
is trivial as well. Since $C_{|\eta_\fq}$ is the restriction
of $h''$ to $|\eta_\fq|$, the claim follows.
\end{pfclaim}

Clearly \eqref{eq_base-change-V} maps each stratum $\sU_\fq$
of the logarithmic stratification of $\sX(\bar\sx)$, to the
corresponding stratum $U_\fq$ of the logarithmic stratification
of $X(\bar x)$ (see \eqref{subsec_log-stratif}).
More precisely, since $\eqref{eq_base-change-V}\times_Y|\eta|$
is ind-\'etale, proposition \ref{prop_generic-boundary}(iii)
implies that the generic point of $\sU_\fq\times_\sY|\eta_\sV|$
gets mapped to the generic point of $U_\fq\times_Y|\eta|$.
We conclude that $E(h)$ contains the generic point of every
stratum $U_\fq\times_Y|\eta|$.

\begin{claim}\label{cl_contains-nuts}
$X(\bar x)\times_Y|\eta|\subset E(h)$.
\end{claim}
\begin{pfclaim} Notice first that
$(X,\underline M)\times_Y|\eta|$ is a regular log scheme
(corollary \ref{cor_smooth-preserve-reg}).

For any geometric point $\xi$ of $X(\bar x)$, denote by
$(X(\xi),\underline M(\xi))$ the strict henselization of
$(X(\bar x),\underline M(\bar x))$ at $\xi$, and set
$C'(\xi):=C'\times_{X(\bar x)}X(\xi)$.
Now, suppose that the support of $\xi$ lies in the stratum
$U_\fq\times_Y|\eta|$, and let $\xi_\fq$ be a geometric
point localized at the generic point of $U_\fq$. By assumption,
$C'(\xi)$ is a finite $X(\xi)$-scheme, tamely ramified along
the non-trivial locus of $(X(\xi),\underline M(\xi))$.
Likewise, $C'(\xi_\fq)$ is tamely ramified along the non-trivial
locus of $(X(\xi_\fq),\underline M(\xi_\fq))$. Pick any
strict specialization map
$(X(\xi_\fq),\underline M(\xi_\fq))\to(X(\xi),\underline M(\xi))$
(see \eqref{subsec_acyclic-log}); it induces a functor
\set\begin{equation}\label{eq_bcov-xi-eta}
\bCov(X(\xi),\underline M(\xi))\to
\bCov(X(\xi_\fq),\underline M(\xi_\fq))
\end{equation}
and theorem \ref{th_Abhyankar_log} implies that $C'(\xi)$
is an object of the source of \eqref{eq_bcov-xi-eta}, which
is mapped, under this functor, to the object $C'(\xi_\fq)$.
By proposition \ref{prop_fund-cuneo}, for any geometric point
$\xi$ of $X(\bar x)\times_Y|\eta|$, the category
$\bCov(X(\xi),\underline M(\xi))$ is equivalent to the category
of finite sets with a continuous action of
$(\underline M{}_\xi)^{\gp\vee}\otimes_\Z\prod_{\ell\neq p}\Z_\ell(1)$.
On the other hand, clearly $\underline M(\bar x)^\sharp$
restricts to a constant sheaf of monoids on $(U_\fq)_\tau$.
In view of \eqref{subsec_punctured}, we deduce that
\eqref{eq_bcov-xi-eta} is an equivalence; lastly, we have
seen that the induced morphism $C'(\xi_\fq)\to X(\xi_\fq)$
is \'etale, {\em i.e.} is a trivial covering, therefore the
same holds for the morphism $C'(\xi)\to X(\xi)$, and
consequently the support of $\xi$ lies in $E(h)$ (claim
\ref{cl_include-pts}). Since $\xi$ is arbitrary, the assertion
follows.
\end{pfclaim}

\begin{claim}\label{cl_make-it-large}
There exists a non-empty open subset $U_Y\subset Y$ such that
$X(\bar x)\times_YU_Y\subset E(h)$.
\end{claim}
\begin{pfclaim} Since $X(\bar x)$ is a noetherian scheme, $E(h)$
is a constructible open subset, hence $Z:=X(\bar x)\!\setminus\!E(h)$
is a constructible closed subset of $X(\bar x)$.
The subset $f_{\bar x}(Z)$ is pro-constructible
(\cite[Ch.IV, Prop.1.9.5(vii)]{EGAIV}) and does not contain
$\eta$ (by claim \ref{cl_contains-nuts}), hence neither
does its topological closure $W$ (\cite[Ch.IV, Th.1.10.1]{EGAIV}).
It is easily seen that $U_Y:=Y\setminus W$ will do.
\end{pfclaim}

\begin{claim}\label{cl_about-Z}
Let $U_Y$ be as in claim \ref{cl_make-it-large}. We have :
\begin{enumerate}
\item
There exists an irreducible closed subset $Z$ of $Y$ of
dimension one, such that
$Z\cap(Y,\underline N)_\tr\cap U_Y\neq\emptyset$.
\item
For any $Z$ as in (i), the induced functor
$\bCov(E(h))\to\bCov(Z\times_YE(h))$ is fully faithful.
\end{enumerate}
\end{claim}
\begin{pfclaim}(i): More generally, let $(A,\fm)$ be any
local noetherian domain of Krull dimension $d\geq 1$, and
$W\subset\Spec\,A$ a proper closed subset; we show that
there exists an irreducible closed subset $Z\subset\Spec\,A$
of dimension one, not contained in $W$. To this aim, we may
assume that $W=\Spec\,A/fA$ for some $f\in\fm\setminus\{0\}$;
let $\fn$ be any maximal ideal of $A[f^{-1}]$, and
$\fp:=A\cap\fn$. Since $A/\fp[f^{-1}]$ is a field,
\cite[Ch.0, Cor.16.3.3]{EGAIV} implies that $Z:=\Spec\,A/\fp$
will do.

(ii): It suffices to check that conditions (i)--(iii) of
proposition \ref{prop_simpler-conds} hold for $Z$ and the
open subset $E(h)$. However, condition (i) is immediate,
since the generic point $\eta_Z$ of $Z$ lies in $U_Y$.
Likewise, condition (ii) holds trivially for the fibre
over the point $\eta_Z$, so it suffices to consider the
fibre over the closed point $y$ of $Z$, in which case the
assertion is just lemma \ref{lem_F_queue}. Lastly, condition
(iii) follows directly from theorem \ref{th_satura-smooth}(iii.b)
and \cite[Ch.IV, Prop.18.8.10, 18.8.12(i)]{EGA4}.
\end{pfclaim}

Let $Z$ be as in claim \ref{cl_about-Z}(i), and endow $Z$
with its reduced subscheme structure. Let also $Z'$ be the
normalization of $Z$; then both $Z$ and $Z'$ are strictly
local, and the morphism $Z'\to Z$ is radicial and surjective,
hence the induced functor
$$
\bCov(Z\times_YE(h))\to\bCov(Z'\times_YE(h))
$$
is an equivalence (lemma \ref{lem_replace}(i)). Taking into
account claim \ref{cl_about-Z}, we are thus reduced to showing
that the morphism
$$
(Z'\times_YE(h))\times_{X(\bar x)}C'\to Z'\times_YE(h)
$$
is a trivial \'etale covering. However, let $\eta_{Z'}$
be the generic point of $Z'$,  and set
$$(Z',\underline N'):=Z'\times_Y(Y,\underline N)
\qquad
(X',\underline M'):=Z'\times_Y(X,\underline M).
$$
The open subset $(Z',\underline N')_\tr$ is dense in $Z'$, by
virtue of claim \ref{cl_about-Z}(i), so the induced
morphism $f':(X',\underline M')\to(Z',\underline N')$ is still
of the type considered in \eqref{sec_acyclic-log}, the 
geometric point $\bar x$ lifts uniquely (up to isomorphism)
to a geometric point $\bar x{}'$ of $X'$, and $h\times_YZ'$
is an object of $\bTame(f',\bar x{}',\kappa(\eta_{Z'}))$
(lemma \ref{lem_tameness-pullsback}(i)).
Hence, we may replace from start $(X,\underline M)$ by
$(X',\underline M')$, $(Y,\underline N)$ by $(Z',\underline N')$,
$h$ by $h\times_YZ'$, after which, we may assume that $Y$
is noetherian and of dimension one. Moreover, taking
into account claim \ref{cl_new-Eprime}, we may assume
that the induced covering $C_{|\eta_\fq}\to|\eta_\fq|$
is trivial, for every maximal point $\eta_\fq$ of
$f^{-1}_{\bar x}(\bar y)$, and it remains to show that
$h$ is trivial under these assumptions.

To this aim, we look at the corresponding commutative diagram
of groups, provided by lemma \ref{lem_horse-with-no-name}(i) :
with the notation of {\em loc.cit.}, we see that in the
current situation, $\beta$ is an isomorphism as well, by
claim \ref{cl_dim-one-ok}, therefore lemma
\ref{lem_horse-with-no-name}(ii) says that the group homomorphism
$\pi_1(|\eta_\fq|_\et,\bar\eta_\fq)\to\pi_1(\bar U_\tr/Y_\et,\xi)$
is surjective, for any maximal point $\eta_\fq$ of
$f^{-1}_{\bar x}(\bar y)$. From this, we deduce that
$\bar h$ is a trivial covering, and therefore there exists
an \'etale covering $C_Y\to|\eta|$ with an isomorphism
$C\isom C_Y\times_{|\eta|}U_\tr$ (\cite[Exp.IX, Th.6.1]{SGA1}).
Denote by $C^\nu_Y$ the normalization of $Y$ in $C_Y$. Also, set
$E':=E(h)\subset\Str(f_{\bar x})$. In light of theorem
\ref{th_satura-smooth}(iii.a) and lemma \ref{lem_F_queue},
it is easily seen that the restriction $E'\to Y$ of $f_{\bar x}$
is surjective; the latter is also a smooth morphism of schemes
(corollary \ref{cor_undercover}(i)). It follows that
$C_Y^\nu\times_YE'$ is the normalization of $E'$ in $C$
(\cite[Ch.IV, Prop.17.5.7]{EGA4}), especially, it is an
\'etale covering of $E'$. We then deduce that $C^\nu_Y$
is already a (trivial) \'etale covering of $Y$
(\cite[Ch.IV, Prop.17.7.1(ii)]{EGA4}), and then clearly $h$
must be a trivial covering as well.
\end{proof}

\begin{remark} Theorem \ref{th_acyclicity-log-smooth} is
the local acyclicity result that gives the name to this
section. However, the title is admittedly not self-explanatory,
and its full justification would require the introduction
of a more advanced theory of the {\em log-\'etale site}, that
lies beyond the bounds of this treatise. In rough terms,
we can try to describe the situation as follows. In lieu
of the standard strict henselization, one should consider
a suitable notion of strict {\em log henselization\/} for
points of the {\em log-\'etale topoi\/} associated with log
schemes. Then, for $f:X\to Y$ as in \eqref{sec_acyclic-log}
with saturated log structures on both $X$ and $Y$, and
log-\'etale points $\tilde x$ of $X$ with image $\tilde y$
in $Y$, one should look, not at our $f_{\bar x}$, but rather
at the induced morphism $f_{\tilde x}$ of strict log
henselizations (of $X$ at $\tilde x$ and of $Y$ at $\tilde y$).
The (suitably defined) {\em log geometric fibres\/} of
$f_{\tilde x}$ will be the {\em log Milnor fibres\/} of $f$ at
the log-\'etale point $\tilde x$, and one can state for such
fibres an acyclicity result : namely, the prime-to-$p$ quotients
of their (again, suitably defined) {\em log fundamental groups\/}
vanish. The proof proceeds by reduction to our theorem
\ref{th_acyclicity-log-smooth}, which, with hindsight, is seen
to supply the essential geometric information encoded in the
more sophisticated log-\'etale language.
\end{remark}

\section{The almost purity toolbox}\label{chap_tool-box}
The sections of this rather heterogeneous chapter are each
devoted to a different subject, and are linked to each other
only very loosely, if they are at all. They have been lumped
here together, because they each contribute a distinct
self-contained little theory, that will find application in
one step or other of the proof of the almost purity theorem
or of its applications in chapter \ref{chap_applications}.
The exception is section \ref{sec_loc-measur-algebras} : it
studies a class of rings more general than the measurable
algebras introduced in section \ref{sec_norm-lengths}; the
results of section \ref{sec_loc-measur-algebras} will not be
used elsewhere in this treatise, but they might be interesting
for other purposes.

Section \ref{subsec_almost-pure-pair} develops the yoga of
almost pure pairs (see definition \ref{def_al-pure}(i)); the
relevance to the almost purity theorem is clear, since the
latter establishes the almost purity of certain pairs $(X,Z)$
consisting of a scheme $X$ and a closed subscheme $Z\subset X$.
This section provides the means to perform various kinds
of reductions in the proof of the almost purity theorem,
allowing to replace the given pair $(X,Z)$ by more tractable
ones.

Section \ref{sec_norm-lengths} introduces measurable (and more
generally, ind-measurable) $K^+$-algebras, for $K^+$ a
fixed valuation ring of rank one : see definitions
\ref{def_working-algeb}(ii) and \ref{def_ind-measur}.
For modules over a measurable algebra, one can define a
well-behaved real-valued normalized length. This length function
is non-negative, and additive for short exact sequences of modules.
Moreover, the length of an almost zero module vanishes (for the
standard almost structure associated with $K^+$). Conversely, under
some suitable assumptions, a module of normalized length zero will
be almost zero.

Lastly, section \ref{sec_fin-groups-quot} studies some questions
concerning the formation of quotients of affine almost schemes
under a finite group action.

\subsection{Non-flat almost structures}\label{sec_non-flat}
This section contains some material that complements the
generalities of \cite[\S2.4, \S2.5]{Ga-Ra} : indeed, whereas
many of the preliminaries in {\em loc.cit.} make no assumptions
on the basic setup $(V,\fm)$ that underlies the whole discussion,
for the more advanced results it is usually required that
$\tilde\fm:=\fm\otimes_V\fm$ is a flat $V$-module.
We shall show how, with more work, one can remove this
condition (or at least, weaken it significantly) and
still recover most of the useful almost ring theory of
\cite{Ga-Ra}. This extension shall be applied in a later
section, in order to state and prove the most general case
of almost purity.

\sset\subsubsection{}\label{subsec_deriv-local}
Let $(V,\fm)$ be any {\em basic setup} in the sense of
\cite[\S2.1.1]{Ga-Ra}, and $R$ any $V$-algebra. For every
interval $I\subset\N$, we have a localization functor
\set\begin{equation}\label{eq_first-on-sC}
\sC^I(R\Mod)\to\sC^I(R^a\Mod)
\qquad
K^\bullet\mapsto K^{\bullet a}
\end{equation}
from complexes of $R$-modules, to complexes of $R^a$-modules,
which is obviously exact, hence it induces a derived localization
functor :
\set\begin{equation}\label{eq_deriv-local}
\sD^I(R\Mod)\to\sD^I(R^a\Mod).
\end{equation}
The functor \eqref{eq_first-on-sC} admits a left (resp. right)
adjoint
\set\begin{equation}\label{eq_left-right}
\sC^I(R^a\Mod)\to\sC^I(R\Mod)
\qquad
K^\bullet\mapsto K^\bullet_!
\qquad
\text{(\ resp. $K^\bullet\mapsto K^\bullet_*$\ )}
\end{equation}
defined by applying termwise to $K^\bullet$ the functor
$M\mapsto M_!$ (resp. $M\mapsto M_*$) for $R^a$-modules
given by \cite[\S2.2.10, \S2.2.21]{Ga-Ra}.
However, if $\tilde\fm$ is not flat, the functor $M\mapsto M_!$
is obviously not exact, and the localization functor
$R\Mod\to R^a\Mod$ does not send injectives to injectives
(cp. \cite[Cor.2.2.24]{Ga-Ra}). This makes it trickier to deal
with constructions in the derived category; for instance, if
$\tilde\fm$ is flat, we get a left adjoint to \eqref{eq_deriv-local},
simply by deriving trivially the exact functor $M\mapsto M_!$.
This fails in the general case, but we shall see later that a
suitable derived version of the construction of $M_!$ is
still available.

\sset\subsubsection{}
For any interval $I\subset\N$, let $\Sigma_I$ be the
multiplicative set of morphisms $\phi$ in $\sD^I(R\Mod)$
such that $\phi^a$ is an isomorphism in $\sD^I(R^a\Mod)$;
arguing as in the proof of \cite[Prop.10.4.4]{We}, one
sees that $\Sigma_I$ is right locally small, hence the
localized category $\Sigma^{-1}_I\sD^I(R\Mod)$ has small
$\Hom$-sets (proposition \ref{prop_calculus-frac}(ii))
which are independent (up to natural isomorphism) of the
choice of universe. Obviously the derived localization
functor factors through a natural functor
\set\begin{equation}\label{eq_Sigma-loc}
\Sigma^{-1}_I\sD^I(R\Mod)\to\sD^I(R^a\Mod).
\end{equation}

\begin{lemma}\label{lem_Sigma-loc}
{\em (i)}\ \
For every interval $I$, the functor \eqref{eq_Sigma-loc}
is an equivalence.
\begin{enumerate}
\addenu
\item
For every $K^\bullet,L^\bullet\in\Ob(\sD^I(R\Mod))$ the induced
$R^a$-linear morphism
$$
\Hom_{\sD^I(R\Mod)}(K^\bullet,L^\bullet)^a\to
\Hom_{\sD^I(R^a\Mod)}(K^{\bullet a},L^{\bullet a})^a
$$
is an isomorphism.
\end{enumerate}
\end{lemma}
\begin{proof}(i): A proof is sketched in \cite[\S2.4.9]{Ga-Ra},
in case $\tilde\fm$ is flat, but in fact this assumption is
superfluous. Indeed, since the unit of adjunction $M\to M_!^a$
is an isomorphism (\cite[Prop.2.2.23(ii)]{Ga-Ra}), it is clear
that the functor $K^\bullet\mapsto K^\bullet_!$ of
\eqref{eq_left-right} descends to a functor
\set\begin{equation}\label{eq_after-Sigma}
\sD^I(R^a\Mod)\to\Sigma^{-1}_I\sD^I(R\Mod)
\end{equation}
such that the composition
$\eqref{eq_Sigma-loc}\circ\eqref{eq_after-Sigma}$ is naturally
isomorphic to the identity automorphism of $\sD^I(R^a\Mod)$.
Likewise, a simple inspection shows that
$\eqref{eq_after-Sigma}\circ\eqref{eq_Sigma-loc}$ is naturally
isomorphic to the identity of $\Sigma^{-1}_I\sD^I(R\Mod)$,
whence the contention.

(ii): The $R$-module $\Hom_{\Sigma_I^{-1}\sD^I(R\Mod)}(K^\bullet,L^\bullet)$
is calculated as the colimit of the system of $R$-modules
$(\Hom_{\sD^I(R\Mod)}(K'^\bullet,L^\bullet)~|~
\phi^\bullet:K'^\bullet\to K^\bullet)$, where $\phi^\bullet$
ranges over the elements of $\Sigma_I$ with target $K^\bullet$.
For such a morphism $\phi^\bullet$, we have
$\fm\cdot H^i(\Cone\,\phi^\bullet)=0$ for every $i\in\Z$.

\begin{claim}\label{cl_depressed}
Let $C^\bullet$ be any object of $\sD^I(R\Mod)$ such
that $(H^iC^\bullet)^a=0$ for every $i\in\Z$. Then
$\Hom_{\sD^I(R\Mod)}(C^\bullet,L^\bullet)^a=0$ for every
$L^\bullet\in\Ob(\sD^I(R\Mod))$.
\end{claim}
\begin{pfclaim} Set
$H_{n,k}:=\Hom_{\sD^I(R\Mod)}(C^\bullet,(\tau^{\geq-n}L^\bullet)[k])$
for every $n,k\in\N$; the natural morphism
$L^\bullet\to\lim_{n\in\N}\tau^{\geq-n}L^\bullet$ is an isomorphism
in $\sC(R\Mod)$, so we have a short exact sequence
$$
0\to\lim_{n\in\N}{}^{\!1}H_{n,-1}\to\Hom_{\sD(R\Mod)}(C^\bullet,L^\bullet)
\to\lim_{n\in\N}H_{n,0}\to 0
$$
(proposition \ref{prop_depressed}(iii)). Hence, it suffices
to check that $H_{n,0}^a=H_{n,-1}^a=0$ for every $n\in\N$, and
we are reduced to the case where $L^\bullet$ is bounded below.

Likewise, set $H'_n:=\Hom_{\sD(R\Mod)}(\tau^{\leq n}C^\bullet,L^\bullet)$
for every $n,k\in\Z$; since the natural morphism
$\colim_{n\in\N}\tau^{\leq n}C^\bullet\to C^\bullet$ is an isomorphism
in $\sC(R\Mod)$, we have an isomorphism
$$
\Hom_{\sD(R\Mod)}(C^\bullet,L^\bullet)\isom
\lim_{n\in\N}H'_n
$$
(proposition \ref{prop_depressed}(ii)), so we may assume
that $C^\bullet$ is a bounded above complex. In this situation,
say that $L^\bullet\in\Ob(\sD^{\geq a}(R\Mod))$ and
$C^\bullet\in\Ob(\sD^{\leq b}(R\Mod))$ for some $a,b\in\Z$;
we may assume that $a\leq b$, and we show, by induction on
$n$, that $(H'_n)^a=0$ for every $n\geq a-1$. The claim will
follow for $n=b$. Indeed, notice that
$$
H'_r=\Hom_{\sD(R\Mod)}(\tau^{\leq r}C^\bullet,\tau^{\leq r}L^\bullet)
\qquad
\text{for every $r\in\Z$}
$$
by virtue of remark \ref{rem_derived-cat}(iv), so the assertion
is clear for $n=a-1$. Next, suppose that the assertion is already
known for every integer $<n$; from the short exact of complexes
$0\to\tau^{\leq n-1}C^\bullet\to\tau^{\leq n}C^\bullet\to(H^nC^\bullet)[-n]\to 0$
we deduce the exact sequence
$$
\Hom_{\sD(R\Mod)}((H^nC^\bullet)[-n],L^\bullet)\to H'_n\to H'_{n-1}
$$
(remark \ref{rem_long-exact-derHom}).
But since $\fm\cdot H^nC^\bullet=0$, the first term in this
exact sequence is annihilated by $\fm$, and the same holds
for the third term, by inductive assumption; then it holds
for the middle term as well, as required.
\end{pfclaim}

By claim \ref{cl_depressed}, we have
$\Hom_{\sD(R\Mod)}(\Cone\,\phi^\bullet,L^\bullet)^a=0$; by
considering the long exact
$\Hom_{\sD(R\Mod)}(-,L^\bullet)$-sequence associated with the
distinguished triangle $K'^\bullet\to K^\bullet\to\Cone\,\phi$,
we deduce an $R^a$-linear isomorphism
$$
\Hom_{\sD^I(R\Mod)}(\phi^\bullet,L^\bullet)^a:
\Hom_{\sD^I(R\Mod)}(K'^\bullet,L^\bullet)^a\isom
\Hom_{\sD^I(R\Mod)}(K^\bullet,L^\bullet)^a
$$
(see remark \ref{rem_long-exact-derHom}). In light of (i),
the assertion follows straightforwardly.
\end{proof}

\sset\subsubsection{}\label{subsec_sketch}
We sketch a few generalities on derived tensor products,
that will be applied to construct useful objects in various
derived categories. Recall first that the tensor product
$-\otimes_R-$ on $R$-modules descends to a bifunctor
$-\otimes_{R^a}-$ (\cite[\S2.2.6, \S2.2.12]{Ga-Ra}), and if
$M$ is a flat $R$-module, then $M^a$ is a flat $R^a$-module
(\cite[Lemma 2.4.7]{Ga-Ra}). It follows that the category
$R^a\Mod$ has enough flat objects, so every bounded above
complex of $R^a$-modules admits a bounded above flat
resolution. Now, given bounded above complexes
$K^\bullet,L^\bullet$ of $R^a$-modules, set
$$
K^\bullet\derotimes_{R^a}L^\bullet:=
(K^\bullet_!\derotimes_RL^\bullet_!)^a
$$
which is a well defined object of $\sD(R^a\Mod)$.

$\bullet$\ \ 
We claim that this rule yields a well defined functor
$$
-\derotimes_{R^a}-:\sD^-(R^a\Mod)\times\sD^-(R^a\Mod)\to
\sD^-(R^a\Mod).
$$
Indeed, suppose $\phi^\bullet:K_1^\bullet\to K^\bullet_2$ is
a morphism in $\sC(R^a\Mod)$, inducing an isomorphism in
$\sD(R^a\Mod)$, and set $C^\bullet:=\Cone\,(\phi^\bullet_!)$;
clearly, $C^{\bullet a}=0$ in $\sD(R^a\Mod)$. Now, pick any
flat bounded above resolution $P^\bullet\to L^\bullet_!$, so
that $K^\bullet_{i!}\otimes_RP^\bullet$ computes
$K^\bullet_{i!}\derotimes_RL^\bullet_!$, for $i=1,2$.
We get natural isomorphism :
$$
\Cone(\phi^\bullet\derotimes_{R^a}L^\bullet)\isom
(C^\bullet\otimes_RP^\bullet)^a\isom
C^{\bullet a}\otimes_{R^a}P^{\bullet a}=0.
$$
so the derived tensor product depends only on the image of
$K^\bullet$ in $\sD^-(R^a\Mod)$; likewise for the argument
$L^\bullet$, whence the contention.

$\bullet$\ \ 
Next, suppose that $P^\bullet\to K^\bullet$ is a bounded
above resolution, with $P^\bullet$ a complex of flat $R^a$-modules;
we claim that there is a natural isomorphism in $\sD(R^a\Mod)$
$$
K^\bullet\derotimes_{R^a}L^\bullet\isom
P^\bullet\otimes_{R^a}L^\bullet.
$$
Indeed, pick any bounded above flat resolution
$Q^\bullet\to L^\bullet_!$; we have natural isomorphisms
$$
K^\bullet\derotimes_{R^a}L^\bullet\isom
(K^\bullet_!\otimes_RQ^\bullet)^a\isom
K^\bullet\otimes_{R^a}Q^{\bullet a}\isom
P^\bullet\otimes_{R^a}Q^{\bullet a}
$$
in $\sD(R^a\Mod)$, where the last holds, since $Q^{\bullet a}$
is a complex of flat $R^a$-modules. Finally, the induced map
$$
P^\bullet\otimes_{R^a}Q^{\bullet a}\to
P^\bullet\otimes_{R^a}L^\bullet
$$
is also an isomorphism in $\sD(R^a\Mod)$, since $P^\bullet$
is a complex of flat $R^a$-modules, so the claim follows.

$\bullet$\ \
Notice that, for any $K^\bullet\in\Ob(\sD^-(R\Mod))$
and any flat resolution $P^\bullet\to K^\bullet$, the induced
morphism $P^{\bullet a}\to K^{\bullet a}$ is a flat resolution;
it follows that, for any $L^\bullet\in\Ob(\sD(R\Mod))$ we get
a natural isomorphism
\set\begin{equation}\label{eq_tensor-and-derive}
(K^\bullet\derotimes_RL^\bullet)^a\isom
K^{\bullet a}\derotimes_{R^a}L^{\bullet a}
\qquad
\text{in $\sD(R^a\Mod)$}.
\end{equation}

\begin{remark}
Clearly, for the derived tensor products of $R^a$-modules,
one has the same commutativity and associativity isomorphisms
as the ones detailed in remark \ref{rem_der-tensor-varie}(i)
for usual modules, as well as the vanishing properties of
lemma \ref{lem_shift-and-shout}(ii).
\end{remark}

We are now ready to return to the question of the existence
of adjoints to derived localization. The key point is the following :

\begin{lemma}\label{lem_derived-tilde}
In the situation of \eqref{subsec_deriv-local}, let $K^\bullet$
be any complex of $R$-modules, $i\in\Z$ any integer and suppose
that :
\begin{enumerate}
\alphaenu
\item
$K^\bullet\in\Ob(\sD^{\leq i}(R\Mod))$
\item
$K^{\bullet a}\in\Ob(\sD^{\leq i-1}(R^a\Mod))$.
\end{enumerate}
Then $\fm\derotimes_VK^\bullet\in\Ob(\sD^{\leq i-1}(R\Mod))$.
\end{lemma}
\begin{proof} We apply the standard spectral sequence
$$
E^2_{pq}:=\Tor^V_p(\fm,H^q K^\bullet)\Rightarrow
H^{q-p}(\fm\derotimes_VK^\bullet).
$$
Indeed, (a) says that $H^qK^\bullet=0$ for every $q>i$, and
(b) says that $(H^iK^\bullet)^a=0$, and therefore
$\fm_V\otimes_VH^iK^\bullet=0$ (\cite[Rem.2.1.4(i)]{Ga-Ra}).
In either case, we conclude that $E^2_{pq}=0$ whenever $q-p\geq i$,
whence the lemma.
\end{proof}

\sset\subsubsection{}\label{subsec_def-fP_i}
Now, let us define inductively :
$$
\fM^\bullet_0:=V[0]
\quad\text{and}\quad
\fM^\bullet_{i+1}:=\fm\derotimes_V\fM^\bullet_i
\quad
\text{for every $i\in\N$}.
$$
A simple induction shows that $\fM^\bullet_i\in\sD^{\leq 0}(V\Mod)$
for every $i\in\N$, so all these derived tensor product are well
defined in $\sD^{\leq 0}(V\Mod)$. Moreover, from the short exact
sequence of $V$-modules
$$
\Sigma
\quad :\quad
0\to\fm\to V\to V/\fm\to 0
$$
we obtain a distiguished triangle
$$
\fM^\bullet_i\derotimes_V\Sigma
\quad :\quad
\fM^\bullet_{i+1}\to\fM^\bullet_i\to
\fM^\bullet_i\derotimes_V(V/\fm)\to\fM^\bullet_{i+1}[1]
\qquad
\text{for every $i\in\N$}.
$$
Especially, we get an inverse system of morphisms in
$\sD^{\leq 0}(V\Mod)$ :
$$
\cdots\to\fM^\bullet_{i+1}\xrightarrow{\ \pi^\bullet_i\ }
\fM^\bullet_i\xrightarrow{\ \pi^\bullet_{i-1}\ }\fM^\bullet_{i-1}
\to\cdots\to\fM^\bullet_0:=V[0].
$$
Also, from \eqref{eq_tensor-and-derive} we deduce natural isomorphisms
in $\sD^{\leq 0}(V^a\Mod)$ :
\set\begin{equation}\label{eq_localize-fM_i}
\fM_i^{\bullet a}\isom V^a[0]
\qquad
\text{for every $i\in\N$}
\end{equation}
and under these identifications, the morphism $\pi_i^{\bullet a}$
corresponds to the identity automorphism of $V^a[0]$ (details
left to the reader). Furthermore, we deduce the following derived
version of \cite[Rem.2.1.4(i)]{Ga-Ra} :

\begin{proposition}\label{prop_derived-tilde}
Let $a,b\in\Z$ be any integers with $a\leq b$. We have :
\begin{enumerate}
\item
For every $K^\bullet\in\Ob(\sD^{[a,b]}(R\Mod))$, the following
conditions are equivalent :
\begin{enumerate}
\item
$K^{\bullet a}\simeq 0$ in $\sD^{[a,b]}(R^a\Mod)$.
\item
$\tau^{\geq a}(\fM^\bullet_{b-a+1}\derotimes_VK^\bullet)\simeq 0$
in $\sD^{[a,b]}(R\Mod)$.
\end{enumerate}
\item
For every morphism $\phi^\bullet:K^\bullet\to L^\bullet$ in
$\sD^{\leq b}(R\Mod)$, the following conditions are equivalent :
\begin{enumerate}
\item
$\phi^{\bullet a}$ is an isomorphism in $\sD^{[a,b]}(R^a\Mod)$.
\item
$\tau^{\geq a}(\fM^\bullet_{b-a+2}\derotimes_V\phi^\bullet)$
is an isomorphism in $\sD^{[a,b]}(R\Mod)$.
\end{enumerate}
\end{enumerate}
\end{proposition}
\begin{proof}(i): From \eqref{eq_localize-fM_i}, it is easily
seen that (b)$\Rightarrow$(a). The other direction follows
straightforwardly from lemma \ref{lem_derived-tilde}, via an
easy descending induction on $b$.

(ii): Again, the direction (b)$\Rightarrow$(a) is immediate
from \eqref{eq_localize-fM_i}. For the other direction, denote
by $C^\bullet$ the cone of $\phi^\bullet$; then
$C^\bullet\in\Ob(\sD^{[a-1,b]}(R\Mod))$, and
$C^{\bullet a}\simeq 0$ in $\sD^{[a-1,b]}(R^a\Mod)$. From (i)
we deduce that
$\tau^{\geq a-1}(\fM^\bullet_{b-a+2}\derotimes_VC^\bullet)\simeq 0$
in $\sD^{[a-1,b]}(R\Mod)$. Since the derived tensor product
is a triangulated functor, the assertion follows easily :
details left to the reader.
\end{proof}

\begin{corollary}\label{cor_derived-tilde}
With the notation of \eqref{subsec_def-fP_i}, the morphism
$$
\tau^{\geq 2-i}\pi^\bullet_i:\tau^{\geq 2-i}\fM^\bullet_{i+1}\to
\tau^{\geq 2-i}\fM^\bullet_i
$$
is an isomorphism in $\sD^{[2-i,0]}(V\Mod)$ for every integer
$i\geq 2$.
\end{corollary}
\begin{proof} By construction, we have a natural isomorphism :
$$
\Cone(\tau^{\geq 2-i}\pi^\bullet_i)\isom
\tau^{\geq 1-i}(\fM^\bullet_i\derotimes_V(V/\fm))
\qquad
\text{in $\sD^{[1-i,0]}(V\Mod)$}
$$
in light of which, the assertion is an immediate consequence
of proposition \ref{prop_derived-tilde}(i).
\end{proof}

\begin{proposition}\label{prop_derive-rl-adjcts}
In the situation of \eqref{subsec_deriv-local}, we have :
\begin{enumerate}
\item
The localization functor $\sD^+(R\Mod)\to\sD^+(R^a\Mod)$
admits the right adjoint :
\set\begin{equation}\label{eq_derived-*}
\sD^+(R^a\Mod)\to\sD^+(R\Mod)
\quad : \quad
K^\bullet\mapsto K^\bullet_{[*]}:=
R\Hom^\bullet_{R^a}(R^a[0],K^\bullet).
\end{equation}
\item
Let $a,b,i\in\Z$ be any three integers with $b\geq a$ and
$i\geq b-a+2$. The localization functor
$\sD^{[a,b]}(R\Mod)\to\sD^{[a,b]}(R^a\Mod)$ admits the left adjoint :
$$
\sD^{[a,b]}(R^a\Mod)\to\sD^{[a,b]}(R\Mod)
\quad :\quad
K^\bullet\mapsto K^\bullet_{[!]}:=
\tau^{\geq a}(\fM^\bullet_i\derotimes_VK^\bullet_{[*]}).
$$
and the right adjoint 
$$
\sD^{[a,b]}(R^a\Mod)\to\sD^{[a,b]}(R\Mod)
\quad :\quad
K^\bullet\mapsto\tau^{\leq b}K^\bullet_{[*]}.
$$
\item
For every $K^\bullet\in\Ob(\sD^+(R^a\Mod))$ (resp.
$L^\bullet\in\Ob(\sD^{[a,b]}(R^a\Mod))$) the counit of adjunction
is an isomorphism
$$
(K^\bullet_{[*]})^a\isom K^\bullet
\qquad
\text{(\ resp.\ $(\tau^{\leq b}L^\bullet_{[*]})^a\isom L^\bullet$ \ ).}
$$
\item
For every $L^\bullet\in\Ob(\sD^{[a,b]}(R^a\Mod))$, the unit of adjunction
is an isomorphism
$$
L^\bullet\to(L^\bullet_{[!]})^a.
$$
\end{enumerate}
\end{proposition}
\begin{proof}(i): This is analogous to lemma \ref{lem_triv-dual}(iii).
Recall the construction : we know that the category $R^a\Mod$ admits
enough injectives (\cite[2.2.18]{Ga-Ra}), hence \eqref{eq_derived-*}
can be represented by $I^\bullet_*$, where $K^\bullet\isom I^\bullet$
is any injective resolution of $K^\bullet$, and $I^\bullet_*$ is
obtained by applying term-wise to $I^\bullet$ the functor $M\mapsto M_*$
of \cite[\S2.2.10]{Ga-Ra}. Indeed, taking into account
\cite[Cor.2.2.19]{Ga-Ra} we get natural isomorphisms :
$$
\begin{aligned}
\Hom_{\sD^+(R\Mod)}(L^\bullet,I^\bullet_*)
\isom\: & H^0\Hom^\bullet_R(L^\bullet,I^\bullet_*) \\
\isom\: & H^0\Hom^\bullet_{R^a}(L^{\bullet a},I^\bullet) \\
\isom\: & \Hom_{\sD(R^a\Mod)}(L^{\bullet a},I^\bullet) \\
\isom\: & \Hom_{\sD(R^a\Mod)}(L^{\bullet a},K^\bullet)
\end{aligned}
$$
for every bounded below complex $L^\bullet$ of $R$-modules.

(iii) follows by direct inspection of the definitions,
taking into account that, for every $R^a$-module $M$, the
counit of adjunction $(M_*)^a\to M$ is an isomorphism
(\cite[Prop.2.2.14(iii)]{Ga-Ra}) : details left to the
reader.

(ii): We consider first the assertion concerning the right
adjoint : let $K^\bullet\in\sD^{[a,b]}(R^a\Mod)$; we may find
an injective resolution $K^\bullet\isom I^\bullet$ such that
$I^j=0$ for every $j<a$, in which case
$I^\bullet_*\in\sD^{\geq a}(R\Mod)$, and $\tau^{\leq b}I^\bullet_*$
represents $\tau^{\leq b}K^\bullet_{[*]}$ in $\sD^{[a,b]}(R\Mod)$.
Then, in view of (i), the assertion is reduced to remark
\ref{rem_derived-cat}(iv). For the left adjoint, let us set
$$
\omega L^\bullet:=\tau^{\geq a}(\fM^\bullet_i\derotimes_VL^\bullet)
\qquad
\text{for every $L^\bullet\in\Ob(\sD^{[a,b]}(R\Mod))$}.
$$
Then $\omega L^\bullet$ is naturally an object of $\sD^{[a,b]}(R\Mod)$,
as explained in remark \ref{rem_der-tensor-varie}(ii), and likewise
for $K^\bullet_{[!]}$, if $K^\bullet$ is any object of
$\sD^{[a,b]}(R^a\Mod)$. We begin with the following :

\begin{claim}\label{cl_about-omegas}
Let $K^\bullet,L^\bullet\in\Ob(\sD^{[a,b]}(R\Mod))$ be any two
objects. Then the natural map
$$
\Hom_{\sD^{[a,b]}(R\Mod)}
(\omega K^\bullet,L^\bullet)\to
\Hom_{\sD^{[a,b]}(R^a\Mod)}(K^{\bullet a},L^{\bullet a})
$$
is an isomorphism.
\end{claim}
\begin{pfclaim} Arguing as in \cite[\S2.2.2]{Ga-Ra}, and taking
into account lemma \ref{lem_Sigma-loc}(i), we reduce to showing
that for any $K^\bullet\in\sD^{[a,b]}(R\Mod)$, the natural morphism
$$
\phi_{K^\bullet}:\omega K^\bullet\to K^\bullet
$$
is initial in the full subcategory of $\sD^{[a,b]}(R\Mod)/K^\bullet$
whose objects are the morphisms $\psi:L^\bullet\to K^\bullet$ that
lie in $\Sigma_{[a,b]}$. However, for any such $\psi$, we have a
commutative diagram in $\sD^{[a,b]}(R\Mod)$ :
$$
\xymatrix{
\omega L^\bullet \ar[r]^-{\phi_{L^\bullet}} \ar[d] &
L^\bullet \ar[d]^\psi \\
\omega K^\bullet \ar[r]^-{\phi_{K^\bullet}} & K^\bullet
}$$
whose left vertical arrow is an isomorphism, by proposition
\ref{prop_derived-tilde}(ii). There follows a morphism
$\phi_{K^\bullet}\to\psi$, and we have to check that this
is the unique morphism from $\phi_{K^\bullet}$ to $\psi$.
However, say that $\alpha,\beta:\phi_{K^\bullet}\to\psi$
are two such morphisms; then their difference is a morphism
$\gamma:=\alpha-\beta:\omega K^\bullet\to L^\bullet$ such
that $\psi\circ\gamma=0$, so $\gamma$ factors through a
morphism $\bar\gamma:\omega K^\bullet\to\Cone\,\psi[-1]$.
Set $C^\bullet:=\tau^{\leq b}\Cone\,\psi[-1]$; then
$C^\bullet\in\sD^{[a,b]}(R\Mod)$, and according to remark
\ref{rem_derived-cat}(iv), $\bar\gamma$ lifts uniquely to
a morphism $\omega K^\bullet\to C^\bullet$ that we denote
again $\bar\gamma$. We deduce a commutative diagram
$$
\xymatrix{
\omega\circ\omega K^\bullet \ar[r] \ar[d]_{\phi_{\omega_{K^\bullet}}}
& \omega C^\bullet \ar[d]^{\phi_{C^\bullet}} \\
\omega K^\bullet \ar[r]^-{\bar\gamma} & C^\bullet.
}$$
Now, by construction $C^{\bullet a}=0$, therefore
$\omega C^\bullet=0$ (proposition \ref{prop_derived-tilde}(i));
on the other hand, $\phi_{\omega_{K^\bullet}}$ is an isomorphism,
by proposition \ref{prop_derived-tilde}(ii). We conclude that
$\bar\gamma=0$, whence $\alpha=\beta$, as sought.
\end{pfclaim}

Assertion (ii.a) is an immediate consequence of (iii) and claim
\ref{cl_about-omegas}; from this, also (iv) is immediate : details
left to the reader.
\end{proof}

\begin{remark}\label{rem_abi}
Let $a,b,i\in\Z$ be any three integers such that $a\leq b$ and
$i\geq b-a+2$. From propositions \ref{prop_derive-rl-adjcts}(ii.a,iii)
and \ref{prop_derived-tilde}(ii) we deduce a natural isomorphism
$$
\tau^{\geq a}(\fM_i^\bullet\derotimes_VK^\bullet)\isom
(K^{\bullet a})_{[!]}
\qquad
\text{for every $K^\bullet\in\Ob(\sD^{[a,b]}(R\Mod))$}
$$
(details left to the reader), which allows to compute
$(K^{\bullet a})_{[!]}$ purely in terms of $K^\bullet$ and
operations within $\sD(R\Mod)$. It turns out that an analogous
isomorphism is available also for $K^{\bullet}_{[*]}$ : this
is contained in the following
\end{remark}

\begin{lemma}\label{lem_*-formula}
Let $a,b,i\in\N$ be any integers such that $a\leq b$ and
$i\geq b-a+2$. For every $K^\bullet\in\Ob(\sD^{[a,b]}(R\Mod))$,
we have a natural isomorphism :
$$
\tau^{\leq b}R\Hom^\bullet_V(\fM_i^\bullet,K^\bullet)
\isom\tau^{\leq b}K^{\bullet a}_{[*]}.
$$
\end{lemma}
\begin{proof} Theorem \ref{th_trivial-dual}(ii) yields a natural
isomorphism
$$
R\Hom_V^\bullet(\fM^\bullet_i,K^\bullet)\isom
R\Hom_R^\bullet(\fM^\bullet_i\derotimes_VR[0],K^\bullet).
$$
To compute the right-hand side, we may fix an injective
resolution $K^\bullet\isom I^\bullet$; the complex $I^{\bullet a}$
is not necessarily injective, but we can find an injective
resolution $\phi:I^{\bullet a}\isom J^\bullet$ (in the category
of bounded below complexes of $R^a$-modules). In view of
\eqref{eq_tensor-and-derive} and \eqref{eq_localize-fM_i}, the
morphism $\phi$ induces a natural transformation
$$
\psi:R\Hom_R^\bullet(\fM^\bullet_i\derotimes_VR[0],K^\bullet)\to
R\Hom^\bullet_{R^a}(R^a[0],K^{\bullet a})=(K^{\bullet a})_{[*]}
$$
and it suffices to show that, for every $j\leq b$, the map
$$
H^j\psi:\Hom_{\sD(R\Mod)}(\fM^\bullet_i\derotimes_VR[0],K^\bullet[j])
\to\Hom_{\sD(R^a\Mod)}(R^a[0],K^{\bullet a}[j])
$$
is an isomorphism. However, by remark \ref{rem_derived-cat}(iv),
the latter is the same as a map
\set\begin{equation}\label{eq_include-truncation}
\Hom_{\sD(R\Mod)}(\fM^\bullet_i\derotimes_VR[0],\tau^{\leq 0}K^\bullet[j])
\to\Hom_{\sD(R^a\Mod)}(R^a[0],\tau^{\leq 0}K^{\bullet a}[j])
\end{equation}
and a direct inspection shows that \eqref{eq_include-truncation}
agrees with the map arising in claim \ref{cl_about-omegas}, for
every $j\leq b$. Especially, for every such $j$, the map
\eqref{eq_include-truncation} is an isomorphism, as sought.
\end{proof}

\begin{proposition}\label{prop_subsume-lr}
Let $a,b,c\in\Z$ be any three integers such that $a\leq b$,
and $K^\bullet,L^\bullet$ any two objects of\/
$\sD^{[a,b]}(R\Mod))$. Suppose that
\begin{enumerate}
\alphaenu
\item
$\Hom_{\sD(R\Mod)}(K^\bullet,X[-j])=0$
for all $j\in[c,b]$ and all $R$-modules $X$ with $X^a=0$.
\item
$\Hom_{\sD(R\Mod)}(Y[-j],L^\bullet)=0$
for all $j\in[a,c]$ and all $R$-modules $Y$ with $Y^a=0$.
\end{enumerate}
Then the natural map
\set\begin{equation}\label{eq_easy-inspect}
\Hom_{\sD(R\Mod)}(K^\bullet,L^\bullet)\to
\Hom_{\sD(R^a\Mod)}(K^{\bullet a},L^{\bullet a})
\end{equation}
is an isomorphism.
\end{proposition}
\begin{proof} We start out with the following observation :

\begin{claim}\label{cl_a'_and-b'}
Consider the following conditions :
\begin{itemize}
\item[(a')]
$\Hom_{\sD(R\Mod)}(K^\bullet,X^\bullet)=0$
for every $X^\bullet\in\Ob(\sD^{[c,b]}(R\Mod))$ such that
$X^{\bullet a}=0$.
\item[(b')]
$\Hom_{\sD(R\Mod)}(Y^\bullet,L^\bullet)=0$
for every $Y^\bullet\in\Ob(\sD^{[a,c]}(R\Mod))$ such that
$Y^{\bullet a}=0$.
\end{itemize}
Then (a)$\Leftrightarrow$(a') and (b)$\Leftrightarrow$(b').
\end{claim}
\begin{pfclaim} Obviously (a')$\Rightarrow$(a). For the
converse, one argues by decreasing induction on $c\leq b$.
Indeed, the case $c=b$ is immediate. Then, suppose that the
sought equivalence has already been established for some
$d\leq b$; if $X^\bullet\in\sD^{[d-1,b]}$ and $X^{\bullet a}=0$,
and if we know that (a) holds with $c:=d-1$, we set
$H^\bullet:=H^cX^\bullet[-c]$, and consider the distinguished
triangle
$$
H^\bullet\to X^\bullet\to\tau^{\geq d}X^\bullet\to H^\bullet[1].
$$
By inductive assumption, we have
$\Hom_{\sD(R\Mod)}(K^\bullet,\tau^{\geq d}X^\bullet)=0$, and
(a) says that
$$
\Hom_{\sD(R\Mod)}(K^\bullet,H^\bullet)=0.
$$
It then follows that $\Hom_{\sD(R\Mod)}(K^\bullet,X^\bullet)=0$,
which shows that the equivalence holds for $c$.

The proof of the equivalence (b)$\Leftrightarrow$(b') is
wholly analogous.
\end{pfclaim}

Fix an integer $i\geq b-a+2$, and set
$$
C^\bullet:=\Cone(
\pi^\bullet_0\circ\cdots\circ\pi^\bullet_i:\fM^\bullet_i\to V[0])
$$
(notation of \eqref{subsec_def-fP_i}); notice that
$C^\bullet\in\sD^{\leq 1}(R\Mod)$, and $C^{\bullet a}=0$.
By virtue of condition (b), claim \ref{cl_a'_and-b'} and
remark \ref{rem_derived-cat}(iv), it follows that :
$$
R^j\Hom_R^\bullet(C^\bullet,L^\bullet)=
\Hom_{\sD(R\Mod)}(\tau^{\geq a}(C^\bullet[-j]),L^\bullet)=0
\qquad
\text{for every $j<c$}
$$
In other words, $D^\bullet:=
R\Hom_R^\bullet(C^\bullet,L^\bullet)\in\sD^{\geq c}(R\Mod)$,
and clearly $D^{\bullet a}=0$; also notice the induced distinguished
triangle
$$
\Sigma
\quad :\quad
D^\bullet\to L^\bullet\to
R\Hom_R^\bullet(\fM^\bullet_i,L^\bullet)\to D^\bullet[-1].
$$
Now, condition (a), claim \ref{cl_a'_and-b'} and remark
\ref{rem_derived-cat}(iv) imply that
$$
\Hom_{\sD(R\Mod)}(K^\bullet,D^\bullet[j])=
\Hom_{\sD(R\Mod)}(K^\bullet,\tau^{\leq b}D^\bullet[j])=0
\qquad
\text{for every $j\leq 0$}
$$
whence, by considering the distinguished triangle
$R\Hom_R^\bullet(K^\bullet,\Sigma)$, natural isomorphisms
$$
\begin{aligned}
\Hom_{\sD(R\Mod)}(K^\bullet,L^\bullet)
\isom\: &
\Hom_{\sD(R\Mod)}(K^\bullet,R\Hom^\bullet_R(\fM^\bullet_i,L)) \\
\isom\: &
\Hom_{\sD(R\Mod)}(\fM^\bullet_i\derotimes_VK^\bullet,L^\bullet) 
& & \text{(by \cite[Th.10.8.7]{We})} \\
\isom\: &
\Hom_{\sD(R\Mod)}
(\tau^{\geq a}\fM^\bullet_i\derotimes_VK^\bullet,L^\bullet) 
& & \text{(by remark \ref{rem_derived-cat}(iv))} \\
\isom\:&
\Hom_{\sD(R^a\Mod)}(K^{\bullet a},L^{\bullet a})
& & \text{(by claim \ref{cl_about-omegas})}
\end{aligned}
$$
whose composition, after a simple inspection, is seen to agree
with the map \eqref{eq_easy-inspect}.
\end{proof}

\begin{remark}\label{rem_a''-and-b''}
(i)\ \
For every interval $I\subset\N$, denote by
$$
\sPhi_I:\sD^I(R/\fm R\Mod)\to\sD^I(R\Mod)
$$
the forgetful functor. It follows easily from lemma
\ref{lem_triv-dual}(iii) and remark \ref{rem_derived-cat}(iv),
that, for every interval $[a,b]$, the functor $\sPhi_{[a,b]}$
admits the right adjoint
$$
\sD^{[a,b]}(R\Mod)\to\sD^{[a,b]}(R/\fm R\Mod)
\qquad
K^\bullet\mapsto\sPsi^r_{[a,b]}K^\bullet:=
\tau^{\leq b}R\Hom_R^\bullet(R/\fm R[0],K^\bullet).
$$

(ii)\ \
Likewise, theorem \ref{th_trivial-dual} and remark
\ref{rem_derived-cat}(iv) imply that $\sPhi_{[a,b]}$
admits the left adjoint
$$
\sD^{[a,b]}(R\Mod)\to\sD^{[a,b]}(R/\fm R\Mod)
\qquad
K^\bullet\mapsto\sPsi^l_{[a,b]}K^\bullet:=
\tau^{\geq a}(K^\bullet\otimes_RR/\fm R[0]).
$$

(iii)\ \
Consider as well the following two conditions :
\begin{itemize}
\item[(a'')]
$\Hom_{\sD(R\Mod)}(K^\bullet,\sPhi_{[c,b]}X^\bullet)=0$
for every $X^\bullet\in\Ob(\sD^{[c,b]}(R/\fm R\Mod))$.
\item[(b'')]
$\Hom_{\sD(R\Mod)}(\sPhi_{[a,c]}Y^\bullet,L^\bullet)=0$
for every $Y^\bullet\in\Ob(\sD^{[a,c]}(R/\fm R\Mod))$.
\end{itemize}
Then, arguing as in the proof of claim \ref{cl_a'_and-b'}
it is easily seen that condition (a) of proposition
\ref{prop_subsume-lr} is equivalent to (a''), and
condition (b) is equivalent to (b'').
\end{remark}

\begin{proposition}\label{prop_ess-image-left}
Let $a,b\in\Z$ be any two integers such that $a\leq b$. For every
object $K^\bullet$ of\/ $\sD^{[a,b]}(R\Mod)$, the following conditions
are equivalent :
\begin{enumerate}
\alphaenu
\item
$K^\bullet$ lies in the essential image of the left adjoint
functor $X^\bullet\mapsto X^\bullet_{[!]}$.
\item
$K^\bullet\derotimes_RR/\fm R[0]\in\Ob(\sD^{<a-1}(R/\fm R\Mod))$.
\end{enumerate}
\end{proposition}
\begin{proof} We start out with the following :

\begin{claim}\label{cl_claim-with-no-name}
We may assume that $V=R$.
\end{claim}
\begin{pfclaim} Clearly condition (b) does not depend
on the underlying ring $V$. It suffices then to remark
that condition (a) depends only on the basic setup
$(R,\fm R)$ (as opposed to the original basic setup
$(V,\fm)$). Indeed, notice that there is a natural equivalence
$$
\Omega:R^a\Mod\isom(R,\fm R)^a\Mod
$$
(where $R^a$ denotes, as in the foregoing, the image of
$R$ in the category of $(V,\fm)$-algebras), and the induced
equivalence of the respective derived categories fits into
an essentially commutative diagram
$$
\xymatrix{ & \sD(R\Mod) \ar[ld] \ar[rd] \\
\sD(R^a\Mod) \ar[rr]^-{\sD(\Omega)} & & \sD((R,\fm R)^a\Mod)
}$$
whose downward arrows are the forgetful functors. Especially,
the left (resp. right) adjoints of these two forgetful functors
share the same essential images.
\end{pfclaim}

Henceforth, we assume that $V=R$ (and therefore, $\fm=\fm R$).
Fix an integer $i\geq b-a+2$, set
$L^\bullet:=\tau^{\geq a-1}(\fM_i^\bullet\derotimes_VK^\bullet)$,
and notice first that, taking into account proposition
\ref{prop_derive-rl-adjcts}(iv) and remark \ref{rem_abi},
condition (a) is equivalent to :
\begin{enumerate}
\addenu\addenu
\alphaenu
\item
The morphism
$\pi^\bullet_0\circ\cdots\circ\pi^\bullet_{i-1}:\fM_i\to V[0]$ induces
an isomorphism $\tau^{\geq a}L^\bullet\to K^\bullet$.
\end{enumerate}

(c)$\Rightarrow$(b): Indeed, set $H:=H^{a-1}L^\bullet$; if (c)
holds, we have a distinguished triangle
$$
H[a-1]\to L^\bullet\to K^\bullet\to H[a]
$$
whence a distinguished triangle in $\sD(V/\fm\Mod)$
\set\begin{equation}\label{eq_L-and-K-and-H}
H[a-1]\derotimes_VV/\fm[0]\to L^\bullet\derotimes_VV/\fm[0]
\to K^\bullet\derotimes_VV/\fm[0]\to H[a]\derotimes_VV/\fm[0].
\end{equation}
However, we have natural isomorphisms
$$
\begin{aligned}
\tau^{\geq a-1}(L^\bullet\derotimes_VV/\fm[0])
\isom\: &
\tau^{\geq a-1}((\fM_i^\bullet\derotimes_VK^\bullet)\derotimes_VV/\fm[0])
& & \text{(by lemma \ref{lem_shift-and-shout}(ii))} \\
\isom\: &
\tau^{\geq a-1}(\fM_i^\bullet\derotimes_V(K^\bullet\derotimes_VV/\fm[0]))
& & \text{(by remark \ref{rem_der-tensor-varie}(i))} \\
\isom\: & \tau^{\geq a-1}(\fM_i^\bullet\derotimes_V
\tau^{\geq a-1}(K^\bullet\derotimes_VV/\fm[0]))
& & \text{(by lemma \ref{lem_shift-and-shout}(ii))} \\
\isom\: & 0
& & \text{(by proposition \ref{prop_derived-tilde}(i))}
\end{aligned}
$$
whence (b), after considering the distinguished triangle
$\tau^{\geq a-1}\eqref{eq_L-and-K-and-H}$.

(b)$\Rightarrow$(a): We remark

\begin{claim}\label{cl_b-and-a''}
Condition (b) is equivalent to condition (a'') of remark
\ref{rem_a''-and-b''}(iii), with $c:=a-1$.
\end{claim}
\begin{pfclaim} Indeed, (b) holds if and only if
$\sPsi^l_{[a-1,b]}K^\bullet=0$ in $\sD^{[a-1,b]}(R/\fm R\Mod)$
(notation of remark \ref{rem_a''-and-b''}(ii)). If the
latter condition holds, then clearly (a'') holds with $c:=a-1$.
Conversely, if (a'') holds for this value of $c$, then
$\Hom_{\sD(R/\fm R\Mod)}(\sPsi^l_{[a-1,b]}K^\bullet,X^\bullet)=0$
for every $X^\bullet\in\sD^{[a-1,b]}(R/\fm R\Mod)$; especially,
the identity automorphism of $\sPsi^l_{[a-1,b]}K^\bullet$
factors through $0$, so $\sPsi^l_{[a-1,b]}K^\bullet$ vanishes.
\end{pfclaim}

From claim \ref{cl_b-and-a''} and remark \ref{rem_a''-and-b''}(iii)
we deduce that if (b) holds, condition (a) of proposition
\ref{prop_subsume-lr} holds for $c:=a-1$, and condition (b)
of the same proposition holds trivially for this value of
$c$, for every $L^\bullet\in\sD^{[a;b]}(R\Mod)$. We conclude
that the natural map
$$
\Hom_{\sD(R\Mod)}(K^\bullet,L^\bullet)\to
\Hom_{\sD(R^a\Mod)}(K^{\bullet a},L^{\bullet a})\isom
\Hom_{\sD(R\Mod)}(K^{\bullet a}_{[!]},L^\bullet)
$$
is an isomorphism, for every $L^\bullet\in\sD^{[a;b]}(R\Mod)$,
whence (a).
\end{proof}

\begin{proposition}\label{prop_ess-image-right}
Let $a,b\in\Z$ be any two integers such that $a\leq b$. For
every $L^\bullet\in\Ob(\sD^{[a,b]}(R\Mod))$, the following
conditions are equivalent :
\begin{enumerate}
\alphaenu
\item
$L^\bullet$ lies in the essential image of the right adjoint
functor $X^\bullet\mapsto\tau^{\leq b}X^\bullet_{[*]}$.
\item
$R\Hom^\bullet_R(R/\fm R[0],L^\bullet)\in\Ob(\sD^{>b+1}(R/\fm R\Mod))$.
\end{enumerate}
\end{proposition}
\begin{proof} By the same argument as in the proof of
claim \ref{cl_claim-with-no-name}, we reduce to the case
where $V=R$. Next, fix $i\in\N$ such that $i\geq b-a+2$, define
$$
K^\bullet:=R\Hom^\bullet_V(\fM_i^\bullet,L^\bullet)
\qquad
\fP_i:=\fM_i\derotimes_VV/\fm[0]
$$
and notice that
\set\begin{equation}\label{eq_conquer-space}
\tau^{\geq a-b-1}\fP_i=0
\qquad
\text{in $\sD^{[a-b-1,0]}(V/\fm\Mod)$}
\end{equation}
due to proposition \ref{prop_derived-tilde}(i). Morever, in view
of proposition \ref{prop_derive-rl-adjcts}(iii) and lemma
\ref{lem_*-formula}, condition (a) is equivalent to :
\begin{enumerate}
\addenu\addenu
\alphaenu
\item
The morphism
$\pi^\bullet_0\circ\cdots\circ\pi^\bullet_{i-1}:\fM_i\to V[0]$
induces an isomorphism $L^\bullet\isom\tau^{\leq b}K^\bullet$.
\end{enumerate}

(c)$\Rightarrow$(b): We argue as in the proof of proposition
\ref{prop_ess-image-left}; namely, set $H:=H^{b+1}K^\bullet$;
if (c) holds, we obtain a distinguished triangle
\set\begin{equation}\label{eq_anew-H-K-L}
H[-b-2]\to L^\bullet\to\tau^{\leq b+1}K^\bullet\to H[-b-1]
\end{equation}
and by considering the induced distinguished triangle
$\tau^{\leq b+1}R\Hom^\bullet_V(V/\fm[0],\eqref{eq_anew-H-K-L})$,
we reduce to observing that
$$
\begin{aligned}
\tau^{\leq b+1}R\Hom^\bullet_V(V/\fm[0],\tau^{\leq b+1}K^\bullet)
\isom\: & \tau^{\leq b+1}
R\Hom_V^\bullet(V/\fm[0],K^\bullet)
& & \text{(by lemma \ref{lem_shift-and-shout})} \\
\isom\: & \tau^{\leq b+1}
R\Hom_V^\bullet(\fP^\bullet_i,L^\bullet)
& & \text{(by \cite[Th.10.8.7]{We})} \\
\isom\: & \tau^{\leq b+1}
R\Hom_V^\bullet(\tau^{\geq a-b-1}\fP^\bullet_i,L^\bullet)
& & \text{(by lemma \ref{lem_shift-and-shout})} \\
\isom\: & 0
& & \text{(by \eqref{eq_conquer-space})}.
\end{aligned}
$$

(b)$\Rightarrow$(a): Again, we proceed as in the proof of the
corresponding assertion in proposition \ref{prop_ess-image-left};
namely, arguing as in the proof of claim \ref{cl_b-and-a''},
we see that condition (b) is equivalent to condition (b'') of
remark \ref{rem_a''-and-b''}(iii), with $c:=b+1$. Hence, if (b)
holds, condition (b) of proposition \ref{prop_subsume-lr} holds
for $c:=b+1$, and notice that condition (a) of {\em loc.cit.}
holds trivially for this value of $c$, for every
$K^\bullet\in\sD^{[a;b]}(R\Mod)$. We conclude that the natural map
$$
\Hom_{\sD(R\Mod)}(K^\bullet,L^\bullet)\to
\Hom_{\sD(R^a\Mod)}(K^{\bullet a},L^{\bullet a})\isom
\Hom_{\sD(R\Mod)}(K^\bullet,\tau^{\leq b}L^{\bullet a}_{[*]})
$$
is an isomorphism, for every $K^\bullet\in\sD^{[a;b]}(R\Mod)$,
whence (a).
\end{proof}

\sset\subsubsection{}\label{subsec_!!}
Let $A$ be any $V^a$-algebra; recall that the localization functor
$V\Alg\to V^a\Alg$ admits a left adjoint $R\mapsto R_{!!}$, whose
restriction to the subcategory of $A_{!!}$-algebras yields a
left adjoint for the localization functor $A_{!!}\Alg\to A\Alg$
(\cite[Prop.2.2.29]{Ga-Ra}). In \cite{Ga-Ra}, we have studied
the deformation theory of $A$-algebras by means of this left
adjoint, under the assumption that $\tilde\fm$ is $V$-flat; here
we wish to show that the same can be repeated in the current
setting, if one appeals instead to the results of the
foregoing paragraphs. To begin with, we remark :

\begin{proposition}\label{prop_back-to-black}
Let $A\to B$ be a morphism of $V^a$-algebras, $N$ a
$B_{!!}$-module. We have:
\begin{enumerate}
\item
If the unit of adjunction $N\to N^a_*$ is injective, the natural map
\set\begin{equation}\label{eq_exalted}
\Exal_{A_{!!}}(B_{!!},N)\to\Exal_A(B,N^a)
\end{equation}
is a bijection (notation of \cite[\S2.5.7]{Ga-Ra}).
\item
If $N^a=0$, then $\Exal_{A_{!!}}(B_{!!},N)=0$.
\end{enumerate} 
\end{proposition}
\begin{proof}(i): Consider any square-zero extension of $A$-algebras
$$
\Sigma
\quad :\quad
0\to N^a\to E\xrightarrow{\ \phi\ } B\to 0.
$$
There follows a square-zero extension of $A_{!!}$-algebras
$$
\Sigma_{!!}
\quad :\quad
0\to N^a_!/\Ker\,\phi_{!!}\to E_{!!}\to B_{!!}\to 0
$$
Under the stated assumption, the counit of adjunction
$N^a_!\to N$ factors uniquely through a $B$-linear map
$g_\phi:N^a_!/\Ker\,\phi_{!!}\to N$; then $g_\phi*\Sigma_{!!}$
(defined as in \cite[\S2.5.5]{Ga-Ra}) yields an element
of $\Exal_{A_{!!}}(B_{!!},N)$ whose image under
\eqref{eq_exalted} equals the class of $\Sigma$. Conversely,
if
$$
\Omega
\quad :\quad
0\to N\to F\xrightarrow{\ \psi\ } B_{!!}\to 0
$$
is a square-zero extension of $A_{!!}$-algebras, then by
adjunction we get a natural map
$$
\Omega^a_{!!}\to\Omega
$$
which in turns, by simple inspection, induces an isomorphism
$g_{\psi^a}*\Omega^a_{!!}\isom\Omega$ in the category of
square-zero $A_{!!}$-algebra extensions of $B_{!!}$ (details
left to the reader). The assertion follows.

(ii): Suppose $N^a=0$, and let $\Omega$ be as in the foregoing;
it follows that $\psi^a:F^a\to B$ is an isomorphism of
$A$-algebras. By adjunction, the morphism $(\psi^a)^{-1}$
corresponds to a map of $A_{!!}$-algebras $\phi:B_{!!}\to E$,
and it is easily seen that $\psi\circ\phi$ is the identity
automorphism of $B_{!!}$, whence the assertion.
\end{proof}

\begin{definition}\label{def_back-to-black}
Let $f:A\to B$ be a morphism of $V^a$-algebras. We set
$$
\L^a_{B/A}:=(\L_{B_{!!}/A_{!!}})^a
$$
which is a simplicial complex of $B$-modules that we call
the {\em almost cotangent complex\/} of $f$.
\end{definition}

\begin{remark}
(i)\ \
In case $\tilde\fm$ is a flat $V$-module, we have introduced in
\cite[Def.2.5.20]{Ga-Ra} a simplicial $B_{!!}$-module $\L_{B/A}$.
Notice that the notation of {\em loc.cit.} agrees with the
current one: indeed, \cite[Prop.8.1.7(ii)]{Ga-Ra} shows that
complex $(\L_{B/A})^a$ obtained by applying the derived localization
functor to $\L_{B/A}$ is naturally isomorphic (in $\sD(s.B\Mod)$)
to the complex of definition \ref{def_back-to-black}.

(ii)\ \
Depending on the context, we will want to regard $\L_{B/A}$ either
as a simplicial object, or as a cochain complex, via the Dold-Kan
isomorphism (\cite[Th.8.4.1]{We}). The resulting slight notational
ambiguity should not be a source of confusion.
\end{remark}

\begin{theorem}\label{th_instead-of-3.2.1}
In the situation of definition {\em\ref{def_back-to-black}},
let $N$ be any $B$-module. Then there are natural isomorphisms
$$
\begin{aligned}
\Der_A(B,N)\isom\: & \Ext^0_B(\L^a_{B/A},N) \\
\Exal_A(B,N)\isom\: & \Ext^1_B(\L^a_{B/A},N).
\end{aligned}
$$
(Notation of\/ \cite[Def.2.5.22(i)]{Ga-Ra}; so, here
we view $\L^a_{B/A}$ as an object of\/ $\sD^{\leq 0}(B\Mod)$.)
\end{theorem}
\begin{proof} The first isomorphism follows easily from
\cite[II.1.2.4.2]{Il} and the natural isomorphism
\set\begin{equation}\label{eq_omegas-are-ok}
\Omega_{B_{!!}/A_{!!}}\isom(\Omega_{B/A})_!
\end{equation}
proved in \cite[Lemma 2.5.29]{Ga-Ra} (the proof in {\em loc.cit.}
does not use the assumption that $\tilde\fm$ is $V$-flat).
Clearly we have
\set\begin{equation}\label{eq_this-is-b}
\Hom_{\sD(B_{!!}\Mod)}(Y[0],N_*[0])=0
\qquad
\text{for every $B_{!!}$-module $Y$ such that $Y^a=0$.}
\end{equation}
On the other hand, we have :

\begin{claim}\label{cl_this-is-a}
$\Hom_{\sD(B_{!!}\Mod)}(\L_{B_{!!}/A_{!!}},X^\bullet)=0$ for every
$X^\bullet\in\Ob(\sD^{[0,1]}(B_{!!}\Mod))$ such that $X^{\bullet a}=0$.
\end{claim}
\begin{pfclaim} For every $X^\bullet\in\Ob(\sD^{[0,1]}(B_{!!}\Mod))$
we have a distinguished triangle
$$
H^0X^\bullet[0]\to X^\bullet\to H^1X^\bullet[-1]\to(H^0X^\bullet)[1]
$$
which reduces to considering the cases where $X^\bullet=M[j]$
for some almost zero $B_{!!}$-module $M$, and $j=0,-1$.
The case where $j=-1$ follows from \cite[III.1.2.3]{Il} and
proposition \ref{prop_back-to-black}(ii). The case where
$j=0$ follows easily from \cite[II.1.2.4.2]{Il} and
\eqref{eq_omegas-are-ok} : details left to the reader.
\end{pfclaim}

Now, \eqref{eq_this-is-b} says that $L^\bullet:=N_*[0]$
fulfills condition (b) of proposition \ref{prop_subsume-lr},
and claim \ref{cl_this-is-a} says that $K^\bullet:=\L_{B_{!!}/A_{!!}}$
fulfills condition (a), so the natural map
$$
\Ext^1_{B_{!!}}(\L_{B_{!!}/A_{!!}},N_*)\to
\Ext^1_B(\L^a_{B/A},N)
$$
is an isomorphism. Taking into account proposition
\ref{prop_back-to-black}(i) and \cite[III.1.2.3]{Il},
the theorem follows.
\end{proof}

\sset\subsubsection{}
For the further study the almost cotangent complex, we shall
need some preliminaries concerning the derived functors of
certain non-additive functors. This material generalizes
the results of \cite[\S8.1]{Ga-Ra}, that were obtained under
the assumption that $\tilde\fm$ is $V$-flat. 

\begin{lemma}\label{lem_simplicial-almost}
Let $(V,\fm)$ be any basic setup, $R$ a simplicial
$V$-algebra, $n\in\N$ an integer, $M$ and $N$ two
$R$-modules such that $H_iM=H_iN=0$ for every
$i\geq n$. The following holds :
\begin{enumerate}
\item
If $M^a=0$ in $\sD(R^a\Mod)$, then $a\cdot\one_M=0$
in $\sD(R\Mod)$, for every $a\in\fm$.
\item
If $\phi:M\to N$ is a morphism of $R$-modules such that
$\phi^a$ in an isomorphism in $\sD(R^a\Mod)$, then for every
$a\in\fm$ we may find a morphism $\psi:N\to M$ in
$\sD(R\Mod)$, such that $\psi\circ\phi=a\cdot\one_M$ and
$\phi\circ\psi=a\cdot\one_N$ in $\sD(R\Mod)$.
\end{enumerate}
\end{lemma}
\begin{proof}(i): For every $R$-module $X$, set
(notation of remark \ref{rem_loop-and-suspend}(iii))
$$
\tau^{\leq-1}X:=\sigma\circ\omega X
$$
According to \cite[I.3.2.1.9(ii)]{Il}, there exists a
natural sequence of morphisms
\set\begin{equation}\label{eq_suspend-loop}
\tau^{\leq-1}X\to X\to s.H_0(X)\to\sigma(\tau^{\leq-1}X)
\qquad
\text{in $\sD(R\Mod)$}
\end{equation}
whose induced sequence of normalized complexes is a
distinguished triangle in $\sD(V\Mod)$ ({\em i.e.} a
distinguished triangle of $\sD(R\Mod)$, in the terminology
of \cite[I.3.2.2.4]{Il}, and in view of \cite[I.3.2.2.5]{Il}).
Let now $n$ and $M$ be as in the lemma; we argue by induction
on $n$. The case where $n=0$ is trivial, so suppose that
$n>0$, and that the assertion has already been proven for
all almost zero $R$-modules $N$ such that $H_iN=0$ for
every $i\geq n-1$. Especially, for $N:=\tau^{\leq 1}X$ and
$P:=s.H_0(M)$ we have $a\cdot\one_{\omega N}=0$ and
$a\cdot\one_P=0$ in $\sD(R\Mod)$, for every $a\in\fm$.
Since the adjunction $\sigma\circ\omega N\to N$ is an
isomorphism in $\sD(R\Mod)$ (\cite[I.3.2.1.10]{Il}), we
deduce that
$$
\Hom_{\sD(R\Mod)}(M,N)^a=0=\Hom_{\sD(R\Mod)}(M,P)^a
$$
whence $\End_{\sD(R\Mod)}(M)^a=0$, by virtue of
\eqref{eq_suspend-loop} (with $X:=M$) and
\cite[I.3.2.2.10]{Il}. The assertion follows.

(ii): Set $C:=\Cone\,\phi$; according to \cite[I.3.2.2]{Il},
we have a distinguished triangle
\set\begin{equation}\label{eq_cone-for-simplices}
M\xrightarrow{\ \phi\ }N\to C\to\sigma M
\qquad
\text{in $\sD(R\Mod)$}
\end{equation}
whence -- by \cite[I.3.2.2.10]{Il} -- an exact sequence
of $V$-modules
$$
\Hom_{\sD(R\Mod)}(N,M)\xrightarrow{\ \alpha\ }
\End_{\sD(R\Mod)}(N)\xrightarrow{\ \beta\ }
\Hom_{\sD(R\Mod)}(N,C).
$$
Now, let us write $a=\sum_{i=1}^na_ib_i$ for some
$a_1,b_1,\dots,a_n,b_n\in\fm$; the assumption on $\phi$
implies that $C^a=0$ in $\sD(R^a\Mod)$, therefore (i) yields
$\beta(a_i\cdot\one_N)=a_i\cdot\beta(\one_N)=0$, so there
exists a morphism $\psi_i:N\to M$ in $\sD(R\Mod)$ such
that $\alpha(\psi_i)=a_i\cdot\one_N$, {\em i.e.}
$\phi\circ\psi_i=a_i\cdot\one_N$ for every $i=1,\dots,n$.
Likewise, by considering the long exact sequence
$\hExt^\bullet_R(\eqref{eq_cone-for-simplices},M)$ provided
by \cite[I.3.2.2.10]{Il} we find, for every $i=1,\dots,n$,
a morphism $\psi'_i:N\to M$ such that
$\psi'_i\circ\phi=b_i\cdot\one_M$. Thus,
$$
a_i\cdot\psi'_i=\psi'_i\circ\phi\circ\psi_i=b_i\cdot\psi_i
\qquad
\text{for every $i=1,\dots,n$}
$$
and a simple computation shows that
$\psi:=\sum_{i=1}^nb_i\cdot\psi_i$ will do.
\end{proof}

\begin{remark}\label{rem_before-prooced}
Before considering non-additive functors, let us see how
to define derived tensor products in $\sD(R^a\Mod)$, for
any simplicial $V$-algebra $R$. We proceed as in
\eqref{subsec_sketch} : for given $R^a$-modules $M,N$, set
$$
M\ellotimes_{R^a}N:=(M_!\ellotimes_RN_!)^a.
$$

(i)\ \
We claim that this rule yields a well defined functor
$$
-\ellotimes_{R^a}-:\sD(R^a\Mod)\times\sD(R^a\Mod)\to\sD(R^a\Mod).
$$
Indeed, say that $\phi:M\to M'$ is a quasi-isomorphism
of $R^a$-modules, and set $C:=\Cone(\phi_!)$. We need to
check that $(\phi_!\ellotimes_RN)^a$ is a quasi-isomorphism,
for any $R$-module $N$; in light of remark
\ref{rem_loop-and-suspend}(iv), it then suffices
to show that $(C\ellotimes_RN)^a=0$; but the latter
$R^a$-module may be computed as
$$
(C\otimes_R\perp^R_\bullet\!N)^a=
C^a\otimes_{R^a}(\perp^R_\bullet\!N)^a
$$
whence the claim, since $C^a=0$ in $R^a\Mod$.

(ii)\ \
Next, suppose that $M$ is a flat $R^a$-module; then
we claim that the natural morphism of $R^a$-modules
$$
M\ellotimes_{R^a}N\to M\otimes_{R^a}N
$$
is a quasi-isomorphism, for every $R^a$-module $N$.
Indeed, by Eilenberg-Zilber's theorem \ref{th_Eilenberg-Zilber},
the assertion comes down to checking that the augmented
simplicial $R^a$-module
$$
(M_!\otimes_R\perp^R_\bullet\!N_!)^a\to M\otimes_{R^a}N
$$
is aspherical. But for every $k\in\N$, the $k$-th column
of the latter is isomorphic to the augmented $R^a$-module
$$
M[k]\otimes_{R^a[k]}(\perp^{R[k]}_\bullet\!N_![k])^a\to
M[k]\otimes_{R^a[k]}N[k]
$$
which is aspherical, since $M[k]$ is a flat $R^a[k]$-module,
whence the assertion.

(iii)\ \
Just as in \eqref{subsec_sketch}, the foregoing immediately
implies yields a natural isomorphism
$$
(M\ellotimes_RN)^a\isom M^a\ellotimes_{R^a}N^a
\qquad
\text{in $\sD(R^a\Mod)$}
$$
for any two $R$-modules $M$ and $N$ (details left to the reader).

(iv)\ \
Furthermore, we get suspension and loop functors
$\sigma$ and $\omega$ for $R^a$-modules, by the rule :
$$
\sigma M:=(\sigma M_!)^a
\quad\text{and}\quad
\omega M:=(\omega M_*)^a
\qquad
\text{for every $R^a$-module $M$}
$$
from which it follows that $\sigma$ is left adjoint to
$\omega$. Then, it is clear that the assertions of remark
\ref{rem_loop-and-suspend}(iii) hold as well for these
functors.

(v)\ \
Likewise, we define the cone of a morphism $\phi:M\to N$
of $R^a$-modules, by the rule
$$
\Cone\,\phi:=(\Cone\,\phi_!)^a
$$
and then the assertion of remark \ref{rem_loop-and-suspend}(iv)
holds also for morphisms of $R^a$-modules.

(vi)\ \
In view of (iv) and (v), it is then easy to check that also
lemma \ref{lem_usual-vanishing} holds {\em verbatim} for
$A:=V$, and any two $R^a$-modules $X$, $Y$.
\end{remark}

\begin{remark}\label{rem-prod-for-dr-a-alg}
(i)\ \ 
In the same vein, we may define derived tensor products
of $R^a$-algebras, for any simplicial $V$-algebra $R$.
Namely, if $S$ and $S''$ are any two $R^a$-algebras, we set
$$
S\ellotimes_{R^a}S':=(S_{!!}\ellotimes_RS'_{!!})^a
$$
(see \eqref{subsec_!!}), where $\ellotimes_R$ denotes the derived
tensor product for $R$-algebras, defined in example
\ref{ex_standard-for-alg}. In view of remarks
\ref{rem_loop-and-suspend}(ii) and \ref{rem_before-prooced}(i),
it is easily seen that this rule defines a functor
$$
-\ellotimes_{R^a}-:\sD(R^a\Alg)\times\sD(R^a\Alg)\to\sD(R^a\Alg)
$$
and moreover, the formation of these tensor products commutes
with the forgetful functor $\sD(R^a\Alg)\to\sD(R^a\Mod)$.

(ii)\ \
Moreover, they are computed by arbitrary flat resolutions :
if $S$ (or $S'$) is a flat $R^a$-algebras, then the natural
morphism
$$
S\ellotimes_{R^a}S'\to S\otimes_{R^a}S'
$$
is an isomorphism in $\sD(R^a\Alg)$. 

(iii)\ \
Furthermore, if $R^a\to S$ is a given morphism of simplicial
$V^a$-algebras, we obtain a well defined functor
$$
\sD(R^a\Alg)\to\sD(S\Alg)
\qquad
S'\mapsto S\ellotimes_{R^a}S'.
$$
Namely, given an $R^a$-algebra $S'$, we pick a resolution
$P\to S'$ with $P$ a flat $R^a$-algebra, and endow
$S\otimes_{R^a}P$ with its natural $S$-algebra structure,
which is independent, up to natural isomorphism, of the
choice of $P$. All the verifications are exercises for the
reader.
\end{remark}

\begin{definition}\label{def_V-homogeneous-fct}
Let $V$ be a ring, $A$ a $V$-algebra, $d\!\in\!\N$ and
$T\!:\!A\Mod\!\to\!A\Mod$ a functor. We say that $T$ is
{\em $V$-homogeneous of degree $d$} if we have
$$
T(a\cdot\one_M)=a^d\cdot\one_{TM}
\qquad
\text{for every $A$-module $M$ and every $a\in V$}.
$$
\end{definition}

\begin{remark}\label{rem_V-homogeneous-funct}
Let $(V,\fm)$ be a basic setup, $A$ a $V$-algebra, $d\in\N$,
and $T:A\Mod\to A\Mod$ a $V$-homogeneous functor of degree $d$.
Suppose that either $d\leq 1$ or else $\fm$ fulfills condition
$(\bB)$ of \cite[\S2.1.6]{Ga-Ra}.

(i)\ \
Let also $\phi:M\to N$ be a homomorphism of $A$-modules such
that $\phi^a:M^a\to N^a$ is an isomorphism (for the almost
structure given by $(V,\fm)$); arguing as in the proof of
lemma \ref{lem_simplicial-almost}(ii), it is easily seen that,
for every $a\in\fm$ there exists a $B$-linear map $\psi:N\to M$
such that $\phi\circ\psi=a\cdot\one_N$ and
$\psi\circ\phi=a\cdot\one_M$. Then, the $V$-homogeneity
property of $T$ implies that $(T\phi)^a$ is an isomorphism
as well (details left to the reader).

(ii)\ \
Hence, $T$ induces a well defined functor
$$
T^a:A^a\Mod\to A^a\Mod
\qquad
M^a\mapsto(TM)^a
\qquad
(\phi^a:M^a\to N^a)\mapsto(T\phi)^a.
$$
For every cardinality $c$, let $\cM_c(A^a)$ be the set of
isomorphism classes of $A^a$-modules which admit a set of
generators of cardinality $\leq c$. We endow $\cM_c(A^a)$
with the uniform structure as in \cite[Def.2.3.1(i)]{Ga-Ra}.
Let $\omega$ be an infinite cardinality, such that $\fm$ is
generated by at most $\omega$ elements; then there exists
a cardinality $\omega'\geq\omega$ such that the isomorphism
class of $T^aM$ lies in $\cM_{\omega'}(A^a)$, for every
$A^a$-module $M$ whose isomorphism class lies in
$\cM_\omega(A^a)$; thus, $T^a$ induces a map
$$
\cM_{\omega,\omega'}(T^a):\cM_\omega(A^a)\to\cM_{\omega'}(A^a)
\qquad
N\mapsto T^aN.
$$
\end{remark}

\begin{lemma}\label{lem_extract-afg-afp}
Let $(V,\fm)$ be a basic setup, $A$ a $V$-algebra, $I,J\subset A$
two ideals, and $\phi:M'\to M$ a morphism of $A^a$-modules with
$I\cdot\Ker\,\phi=J\cdot\Coker\,\phi=0$. Then there exists an
$A^a$-linear morphism $\lambda:I\otimes_AJ\otimes_AM\to M'$ with
$$
\phi\circ\lambda=\mu_{I,J}\otimes_A\one_M
\qquad\text{and}\qquad
\lambda\circ(I\otimes_AJ\otimes_A\phi)=\mu_{I,J}\otimes_A\one_{M'}
$$
where $\mu_{I,J}:I\otimes_AJ\to A$ is the multiplication law :
$a\otimes b\mapsto ab$ for every $a\in I$ and $b\in J$.
\end{lemma}
\begin{proof} Let
$M'\xrightarrow{\bar\phi}M_0:=\Img\,\phi\xrightarrow{i}M$ be
the natural factorization of $\phi$; then for every $a\in I$
and $b\in J$ we have $A^a$-linear morphisms
$$
\psi_a:M_0\to M'
\quad\text{and}\quad
\mu_b:M\to M_0
\qquad\text{such that}\qquad
\psi_a\circ\bar\phi=a\cdot\one_{M'}
\quad\text{and}\quad
i\circ\mu_b=b\cdot\one_M.
$$
Let us check that $\bar\phi\circ\psi_a=a\cdot\one_{M_0}$; since
$\bar\phi$ is an epimorphism, it suffices to show that
$\bar\phi\circ\psi_a\circ\bar\phi=(a\cdot\one_{M_0})\circ\bar\phi$,
which is clear. Likewise, let us check that
$\mu_b\circ i=b\cdot\one_{M_0}$; since $i$ is a monomorphism, it
suffices to show that $i\circ\mu_b\circ i=i\circ b\cdot\one_{M_0}$,
which is clear. For every $a\in I$ and $b\in J$ set
$\lambda_{a,b}:=\psi_a\circ\mu_b:M\to M'$; we deduce that
$$
\phi\circ\lambda_{a,b}=i\circ\bar\phi\circ\psi_a\circ\mu_b=
i\circ(a\cdot\one_{M_0})\circ\mu_b=i\circ\mu_b\circ(a\cdot\one_M)=
ab\cdot\one_M
$$
and a similar computation yields :
$\lambda_{a,b}\circ\phi=ab\cdot\one_{M'}$. Thus, we obtain a map
$$
I\times J\times M\to M'
\qquad
(a,b,x)\mapsto\lambda_{a,b}(x).
$$
A simple inspection shows that this map is $A$-trilinear, hence
it factors uniquely through an $A$-linear map
$\lambda:I\otimes_AJ\otimes_AM\to M'$, and the foregoing easily
implies that $\lambda$ fulfills the required identities.
\end{proof}

\begin{proposition}\label{prop_ext-std-to-m-nonflat}
In the situation of remark {\em\ref{rem_V-homogeneous-funct}} the
following holds :
\begin{enumerate}
\item
If $T$ commutes with filtered colimits, and for every free
$A$-module $L$ of finite rank, $(TL)^a$ is a flat $A^a$-module
(resp. is the zero $A^a$-module), then for every flat $A^a$-module
$M$, the $A^a$-module $T^aM$ is flat (resp. is the zero $A^a$-module).
\item
If for every free $A$-module $L$ of finite rank (resp. for every
free $A$-module $L$), $(TL)^a$ is an almost projective $A^a$-module,
then for every almost projective and almost finitely generated
$A^a$-module $M$ (resp. for every almost projective $A^a$-module
$M$), the $A^a$-module $T^aM$ is almost projective.
\item
The map $\cM_{\omega,\omega'}(T^a)$ of remark
{\em\ref{rem_V-homogeneous-funct}(ii)} is uniformly continuous.
\end{enumerate}
\end{proposition}
\begin{proof}(i): Say that $M=N^a$ for an $A$-module $N$, and
let $\cF$ be the category whose objects are the pairs $(C,\phi)$,
where $C$ is a finitely presented $A$-module and $\phi:C\to N$
is an $A$-linear map; the morphisms $\psi:(C,\phi)\to(C',\phi')$
in $\cF$ are the $A$-linear maps $\psi:C\to C'$ such that
$\phi'\circ\psi=\phi$, with the obvious composition law. We have
a functor
$$
F:\cF\to A\Mod
\qquad
(C,\phi)\mapsto C
\qquad
((C,\phi)\xrightarrow{\psi}(C',\phi'))\mapsto(C\xrightarrow{\psi}C')
$$
as well as a natural co-cone $\beta:F\Rightarrow c_N$ given
by the rule : $(C,\phi)\mapsto\phi$. It is easily seen that
$\cF$ is filtered, and $\beta$ is a universal co-cone : the
details are left to the reader. Since $T$ commutes with
filtered colimits, $TN$ is isomorphic to the colimit of the
functor $T\circ F:\cF\to A\Mod$, and
$T*\beta:T\circ F\Rightarrow c_{TN}$ is again a universal
co-cone. We need to show that $\fm\Tor^A_1(TN,X)=0$ for every
$A$-module $X$ (resp. that $\fm TN=0$), and by the foregoing
it then suffices to check that for every $(C,\phi)\in\Ob(\cF)$,
the induced map
$$
\lambda:=\Tor_1^A(T\phi,X):\Tor_1^A(TC,X)\to\Tor_1^A(TN,X)
$$
is almost zero (resp. that $T\phi$ is almost zero). However,
since $N^a$ is a flat $A^a$-module, \cite[Lemma 2.4.17]{Ga-Ra}
implies that for every $a\in\fm$ there exists a free $A$-module
$L$ of finite rank and $A$-linear maps $\phi':C\to L$,
$\phi'':L\to M$ with $a\cdot\phi=\phi''\circ\phi'$. Hence
$a^d\cdot T\phi=T(\phi'')\circ T(\phi'')$, and therefore
$a^d\cdot\lambda^a=0$, since by assumption $(TL)^a$ is a flat
$A^a$-module (resp. and therefore $a^d\cdot(T\phi)^a=0$, since
$(TL)^a=0$); the contention follows, since we have either $d\leq 1$,
or else $\fm$ fulfills condition $(\bB)$.

(ii): Suppose first that for every free $A$-module $L$ of
finite rank, $(TL)^a$ is almost projective, and let $M$
be an almost projective and almost finitely generated $A$-module.
According to \cite[2.4.15]{Ga-Ra}, for every $a\in\fm$ there
exists a free $A$-module $L$ of finite rank and $A^a$-linear
morphisms $\phi':M\to L^a$, $\phi'':L^a\to M$ such that
$a\cdot\one_M=\phi''\circ\phi'$. We deduce that
$$
a^d\cdot\Ext^1_A(\one_{TM},X)^a=
\Ext^1_A(T\phi',X)^a\circ\Ext^1_A(T\phi'',X)^a=0
$$
for every $A$-module $X$, since $\Ext^1_A(TL,X)^a=0$ by assumption
(and by \cite[Rem.2.4.12(i)]{Ga-Ra}). The contention follows, again
because either $d\leq 1$ or else $\fm$ fulfills condition $(\bB)$. One
argues similarly in case $M$ is almost projective and $(TL)^a$ is
almost projective for every free $A$-module $L$ : the details shall
be left to the reader.

(iii): Let $\fm_0\subset\fm$ be a subideal of finite type,
and $\phi:M'\to M$ a morphism of $A^a$-module with
$\fm_0\Ker\,\phi=\fm_0\Coker\,\phi=0$. Denote by
$\fm_0^{(2d)}\subset\fm_0$ the subideal generated by the system
$(a^d~|~a\in\fm^2_0)$. By lemma \ref{lem_extract-afg-afp}, for
every $a\in\fm^2_0$ there exists a morphism $\lambda:M\to M'$
of $A^a$-modules with $\phi\circ\lambda=a\cdot\one_M$ and
$\lambda\circ\phi=a\cdot\one_{M'}$; therefore
$$
T^a(\phi)\circ T^a(\lambda)=a^d\cdot\one_{T^aM}
\qquad\text{and}\qquad
T^a(\lambda_a)\circ T^a(\phi)=a^d\cdot\one_{T^aM'}
$$
which implies that $\fm_0^{(2d)}\cdot\Ker\,T^a(\phi)=
\fm_0^{(2d)}\cdot\Coker\,T^a(\phi)=0$, whence the assertion, since
we either have $d\leq 1$ or $\fm$ fulfills condition $(\bB)$.
\end{proof}

\begin{example} Let $(V,\fm)$ be a basic setup such that $\fm$
satisfies condition $(\bB)$; let $A$ be a $V$-algebra, and $P$
any $A$-module. According to \cite[Ch.X, \S9, n.3]{BouAH}, for
every integer $n>0$ we have a complex of $A$-modules
$$
\Sigma^\bullet_A(P)
\qquad : \qquad
0\to\Lambda^n_AP\to\cdots\to(\Lambda_A^jP)\otimes_A(\Sym^{n-j}_AP)
\to\cdots\to\Sym^n_AP\to 0
$$
which is functorial for morphisms of $A$-modules, and is acyclic
if $P$ is a flat $A$-module (\cite[Ch.X, \S9, n.3, Prop.3]{BouAH}).
Clearly, for every $j\in\N$ the functor
$$
T^j:A\Mod\to A\Mod
\qquad
P\mapsto H^j(\Sigma^\bullet_A(P))
$$
fulfills the conditions of \eqref{subsec_homogeneous-functor}
(with $d:=n$) and $T^jP=0$ if $P$ is a flat $A$-module.
By proposition \ref{prop_ext-std-to-m-nonflat}(i), we deduce
that the induced complex of $A^a$-modules $\Sigma^\bullet_A(P)^a$
is acyclic, if $P^a$ is a flat $A^a$-module. With this observation,
all the results of \cite[\S4.3, \S4.4]{Ga-Ra} extend now
{\em verbatim} to the case where $\fm$ fulfills condition
$(\bB)$ (whereas the proofs in {\em loc.cit.} used the stronger
condition that $\fm$ is a flat $V$-module : the details shall
be left to the reader).
\end{example}

\sset\subsubsection{}\label{subsec_homogeneous-functor}
Let $\pi:V\tdu\AlgMod\to V\Alg$ be the functor such that
$\pi(A,M):=A$ and $\pi(f,\phi):=f$ for every object $(A,M)$
and every morphism $(f,\phi)$ of $V\tdu\AlgMod$. We consider
a functor
$$
T:V\tdu\AlgMod\to V\tdu\AlgMod
\qquad
\text{such that $\pi\circ T=\pi$}.
$$
Let also be $R$ a simplicial $V$-algebra, $d\in\N$ and suppose that :
\begin{itemize}
\item
$T$ is {\em $V$-homogeneous of degree $d$, i.e.} for every
$V$-algebra $A$, the restriction of $T$
$$
T_A:A\Mod\to A\Mod
$$
is $V$-homogeneous of degree $d$.
\item
$T_A$ commutes with filtered colimits (cp. \eqref{subsec_new-subsection})
for every $V$-algebra $A$.
\item
Either $d\leq 1$ or else the ideal $\fm\cdot H_0(R)$ of $H_0(R)$
satisfies condition $(\bB)$ of \cite[\S2.1.6]{Ga-Ra}.
\end{itemize}

\begin{remark}\label{rem_homogeneous_functor}
(i)\ \
As detailed in \eqref{subsec_new-subsection}, the functor $T$
induces a functor
$$
T_R:R\Mod\to R\Mod
\qquad
(M[n]~|~n\in\N)\mapsto(T_{R[n]}M[n]~|~n\in\N).
$$

(ii)\ \
Let $T$ and $T'$ be two functors as in
\eqref{subsec_homogeneous-functor}, $V$-homogeneous of degrees
respectively $d$ and $d'$; then we get the functor $V$-homogeneous
of degree $d+d'$
$$
T\otimes T':V\tdu\AlgMod\to V\tdu\AlgMod
\qquad
(A,M)\mapsto(A,T_AM\otimes_AT'_AM)
$$
and it is easily seen that $(T\otimes T')_A$ commutes
with filtered colimits, for every $V$-algebra $A$.

(iii)\ \
Likewise, if $f:T\to T'$ is a natural transformation
of functors fulfilling the conditions of
\eqref{subsec_homogeneous-functor}, for the same degree $d$;
then the same holds for the functors $\Ker\,f$ and $\Coker\,f$.

(iv)\ \
In the situation of \eqref{subsec_homogeneous-functor}, let
$f:H_0(R)\to B$ be any morphism of $V$-algebras. If condition
$(\bB)$ holds for $\fm\cdot H_0(R)$, it holds also also for
$\fm B$; then remark \ref{rem_V-homogeneous-funct}(ii) implies
that $T_B$ induces a functor $T_{B^a}:=T^a_B:B^a\Mod\to B^a\Mod$.
Especially, the above holds with $B:=R[p]$, for any $p\in\N$,
and $f$ the natural map given by the degeneracies of the
simplicial $V$-algebra $R$. It is then clear that $T_R$ induces
a well defined functor
$$
T_{R^a}:R^a\Mod\to R^a\Mod.
$$
The following result shows that, in this situation, the
construction of left derived functors descends likewise
to $R^a$-modules.
\end{remark}

\begin{theorem}\label{th_pluri-vanish}
In the situation of \eqref{subsec_homogeneous-functor},
the following holds :
\begin{enumerate}
\item
Let $\phi:M\to N$ be a morphism of $R$-modules, and $n\in\N$
an integer such that $H_i(\phi)^a:(H_iM)^a\to(H_iN)^a$
is an isomorphism of $V^a$-modules, for every $i\leq n$. Then
$$
H_i(LT\phi)^a:H_i(LTM)^a\to H_i(LTN)^a
$$
is an isomorphism of $V^a$-modules, for every $i\leq n$.
\item
Especially, the functor $T_R$ induces a well defined left
derived functor
$$
LT_{R^a}:\sD(R^a\Mod)\to\sD(R^a\Mod).
$$
\end{enumerate}
\end{theorem}
\begin{proof} Clearly (ii) follows from (i).

(i): In light of corollary \ref{cor_truncate-and-derive}(i),
we may replace $M$ and $N$ by respectively $\cosk_nM$ and
$\cosk_nN$, after which, we may assume that $H_iM=H_iN=0$
for every $i>n$, so $\phi^a$ is an isomorphism in $\sD(R^a\Mod)$.
Next, by proposition \ref{prop_derive-non-add}, we may replace
$M$ and $N$ by their respective standard free resolutions,
in which case we are reduced to checking that the induced
map $(T_R\phi)^a$ is an isomorphism in $\sD(R^a\Mod)$. Since
$T$ is homogeneous and $(\bB)$ holds for $\fm\cdot H_0(R)$,
the latter assertion follows straightforwardly from lemma
\ref{lem_simplicial-almost}(ii).
\end{proof}

\sset\subsubsection{}\label{subsec_drop_B}
Next, we wish to decide how much of the foregoing theory
can be salvaged, when we drop condition ($\bB$). We will
concentrate on the functors that are relevant to the later
study of the cotangent complex. Thus, henceforth, for
every integer $d\in\N$, we shall denote by $T^d$ one of
the three standard functors
$$
\Sym^d,\Lambda^d,\Gamma^d:V\tdu\AlgMod\to V\tdu\AlgMod
$$
(namely, the symmetric and exterior $d$-th power functors,
and the $d$-th divided power functor). Notice that the
functor $T^d$ fulfills the conditions of
\eqref{subsec_homogeneous-functor}, for every $d\in\N$.
Moreover, in all three cases we have natural identifications :
\set\begin{equation}\label{eq_T_1}
T^1\isom\one_{V\tdu\AlgMod}.
\end{equation}
In general, we can no longer expect that $T^d$ descends
to almost modules; indeed, consider the following :

\begin{example}\label{ex_worst}
Let $(V,\fm)$ be as in \eqref{subsec_drop_B}, and $p\in\N$
any prime integer. The Frobenius map $\Phi_R:R/pR\to R/pR$
for $V$-algebras $R$, can be seen as a homogeneous polynomial
law of degree $p$ on the $V$-module $V/pV$. Denote by
$\fm^{(p)}\subset V$ the ideal generated by $(x^p~|~x\in\fm)$.
Clearly $\Phi$ descends to a polynomial law
$$
\bar\Phi:V/(pV+\fm)\to W_p:=V/(pV+\fm^{(p)})
$$
which is still homogeneous of degree $p$, so it factors
through a unique $V$-linear map
$$
\Gamma^p_A(V/(pV+\fm))\to W_p.
$$
The latter is surjective, since its image contains the class
of the unit element of $V$. Now, the proof of
\cite[Prop.2.1.7(ii)]{Ga-Ra} shows that -- if $(\bB)$ does
not hold for $\fm$ -- there exists some prime $p$ such that
$pV+\fm^{(p)}$ does not contain $\fm$, therefore $W_p^a\neq 0$.
This shows that the functor
$$
V\Mod\to V^a\Mod
\qquad
M\mapsto(\Gamma^p_AM)^a
$$
does not factor through $V^a\Mod$.
\end{example}

However, we will see that example \ref{ex_worst} is, in
a sense, the worst that can happen.

\begin{lemma}\label{lem_i-j-SGL}
In the situation of \eqref{subsec_drop_B}, we have :
\begin{enumerate}
\item
There are natural transformations
$$
T^{i+j}\to T^i\otimes T^j\to T^{i+j}
\qquad
\text{for every $i,j\in\N$}
$$
whose composition equals $\binom{i+j}{i}\cdot\one_{T^{i+j}}$.
(Notation of remark {\em\ref{rem_homogeneous_functor}(ii)}.)
\item
For every simplicial $V$-algebra $R$, the maps of {\em(i)}
induce natural transformations
$$
LT_R^{i+j}\to LT_R^i\ellotimes_R LT_R^j\to LT_R^{i+j}
\qquad
\text{for every $i,j\in\N$}
$$
whose composition equals $\binom{i+j}{i}\cdot\one_{LT_R^{i+j}}$.
\end{enumerate}
\end{lemma}
\begin{proof}(i): For $T=\Lambda$, the sought morphism
$\Lambda^{i+j}\to\Lambda^i\otimes\Lambda^j$ is the one
denoted $\Delta_{i,j}$ in \cite[\S4.3.20]{Ga-Ra}, and
the morphism $\Lambda^i\otimes\Lambda^j\to\Lambda^{i+j}$
is given by the rule : $x\otimes y\mapsto x\wedge y$,
for every $(A,M)\in\Ob(V\tdu\AlgMod)$, every
$x\in\Lambda^i_AM$ and every $y\in\Lambda^j_AM$. The
sought identity follows easily from the explicit
formula given in \cite[(4.3.21)]{Ga-Ra} : details left
to the reader.

For $T=\Sym$, the natural transformation
$\Sym^{i+j}\to\Sym^i\otimes\Sym^j$ is given by a similar
formula; namely, for every object $(A,M)$ of $V\tdu\AlgMod$,
every sequence of elements $x_1,\dots,x_{i+j}\in M$, and
every subset $I\subset\{1,\dots,i+j\}$, set
$x_I:=\prod_{k\in I}x_k\in\Sym^k_AM$ (where the
multiplication is formed in the graded ring
$\Sym^\bullet_AM$); one checks easily that the rule
$$
x_1\cdots x_{i+j}\mapsto\sum_{I,J}x_I\otimes x_J
$$
defines a well defined $A$-linear map on $\Sym^{i+j}_AM$
(where the sum ranges over all the partitions $(I,J)$ of
$\{1,\dots,i+j\}$ such that the cardinality of $I$
equals $i$). Then one defines a map
$\Sym^i_AM\otimes_A\Sym^j_AM\to\Sym^{i+j}_AM$ by the
rule : $u\otimes v\mapsto u\cdot v$ for every
$u\in\Sym^i_AM$ and $v\in\Sym^j_AM$. Again, the sought
identity is verified by direct computation.

If $T=\Gamma$, then for every object $(A,M)$
of $V\tdu\AlgMod$ we define a homogeneous polynomial law
$M\squig\Gamma^i_AM\otimes_A\Gamma^j_AM$ of degree $i+j$,
by the rule : $x\mapsto x^{[i]}\otimes x^{[j]}$ for every
$A$-algebra $B$, and every $x\in B\otimes_AM$.
This law yields a transformation
$\Gamma^{i+j}\to\Gamma^i\otimes\Gamma^j$ as sought. Next,
the multiplication of the graded ring functor $\Gamma^\bullet$
gives a transformation $\Gamma^i\otimes\Gamma^j\to\Gamma^{i+j}$.
The composition of these two transformations is characterized
as the unique transformation $\psi$ of $\Gamma^{i+j}$ such that
$\psi_{(A,M)}(x^{[i+j]})=x^{[i]}\cdot x^{[j]}$ for every
object $(A,M)$ of $V\tdu\AlgMod$, and every $x\in M$.
Then, by corollary \ref{cor_divided-powers} we must have
$\psi_{(A,M)}=\binom{i+j}{i}\one_M$, as required.

(ii): Recall that the three functors $\Lambda$, $\Gamma$
and $\Sym$ transform free $A$-modules into free $A$-modules,
for every $V$-algebra $A$. Therefore, the sought natural
transformation is none else than the map
$$
T^{i+j}_R(\perp^R_\bullet\!M)\to
T^i_R(\perp^R_\bullet\!M)\otimes_RT^j_R(\perp^R_\bullet\!M)\to
T^{i+j}_R(\perp^R_\bullet\!M)
$$
given by (i), for any $R$-module $M$.
\end{proof}

\begin{proposition}\label{prop_small-degrees}
In the situation of \eqref{subsec_drop_B}, let $p\in\N$ be a
prime integer, $R$ a simplicial $V\otimes_\Z\Z_{(p)}$-algebra,
$M$ an $R$-module such that $M^a=0$ in $\sD(R^a\Mod)$.
Then
$$
(LT^d_RM)^a=0
\qquad
\text{in $\sD(R^a\Mod)$, for every $d\in\N$ such that $(p,d)=1$}.
$$
\end{proposition}
\begin{proof} In light of corollary \ref{cor_truncate-and-derive}(i),
we may replace $M$ by $\cosk_iM$, in which case lemma
\ref{lem_simplicial-almost}(i) implies that $a\cdot\one_M=0$
for every $a\in\fm$. Now, we apply lemma \ref{lem_i-j-SGL}(ii)
with $i=1$ and $j=d-1$ (notice that $j\geq 0$, since $d>0$);
in view of \eqref{eq_T_1}, there result natural transformations
$$
LT^d_RM\to M\otimes_RLT^{d-1}_RM\to LT^d_RM
$$
whose composition is $d\cdot\one_{LT^d_RM}$. Since the image
of $d$ is invertible in $R$, it follows easily that
$a\cdot\one_{LT^d_RM}=0$ for every $a\in\fm$, whence the
assertion.
\end{proof}

\begin{theorem}\label{th_prime-vanishing}
Let $d,n,p\in\N$ be any three integers, with $d>0$ and
$p$ a prime. Let also $R$ be a simplicial
$V\otimes_\Z\Z_{(p)}$-algebra, $M$ an $R$-module,  and
suppose that
\begin{enumerate}
\alphaenu
\item
$H_iM=0$ for every $i<n$.
\item
$M^a=0$ in $\sD(R^a\Mod)$.
\end{enumerate}
Then we have :
\begin{enumerate}
\item
$H_i(L\Gamma^d_RM)^a=0$ for every $i<n$.
\item
$H_i(L\Lambda^d_RM)^a=0$ for every $i<n+p-1$.
\item
$H_i(L\Sym^d_RM)^a=0$ for every $i<n+2(p-1)$.
\end{enumerate}
\end{theorem}
\begin{proof}(i) is just a special case of corollary
\ref{cor_truncate-and-derive}(ii), and for $d<p$,
assertions (ii) and (iii) follow from the more
general proposition \ref{prop_small-degrees}, hence we
may assume that $d\geq p$. For every $j=0,\dots,d$, set
$$
F^j:=\Gamma^j\otimes\Lambda^{d-j}
\qquad
G^j:=\Lambda^j\otimes\Sym^{d-j}.
$$
By remark \ref{rem_homogeneous_functor}(ii), these functors
are homogeneous of degree $d$, and fulfill the conditions
of \eqref{subsec_homogeneous-functor}. Moreover,
since $\Gamma$, $\Lambda$ and $\Sym$ transform free
modules into free modules, we have natural isomorphisms
of functors :
\set\begin{equation}\label{eq_decompose_F}
LF^j_R\isom L\Gamma^j_R\ellotimes_RL\Lambda^{d-j}_R
\quad
LG^j_R\isom L\Lambda^j_R\ellotimes_RL\Sym^{d-j}_R
\qquad
\text{for every $j=0,\dots,d$}
\end{equation}
(details left to the reader).

(ii): According to \cite[I.4.3.1.7]{Il}, there is a natural
complex of functors :
\set\begin{equation}\label{eq_cplx-of-functs}
0\to F^d_R\xrightarrow{\ \partial_{d-1}}F_R^{d-1}\to\cdots\to
F_R^1\xrightarrow{\ \partial_0}F_R^0\to 0
\end{equation}
which is exact on flat $R$-modules. Set
$Z_p:=\Ker\,\partial_{p-1}$; due to remark
\ref{rem_homogeneous_functor}(iii) we may consider the derived
functor $LZ_p$, and from corollary
\ref{cor_truncate-and-derive}(ii) and assumption (a), we get
\set\begin{equation}\label{eq_boundary-functor}
H_i(LZ_pM)=0
\qquad
\text{for every $i<n$}.
\end{equation}
The evaluation of \eqref{eq_cplx-of-functs} on
$(\perp^R_\bullet\!M)^\Delta$ gives an exact sequence of $R$-modules :
$$
0\to LZ_pM\to LF_R^{p-1}M\to\cdots\to LF^1_RM\to LF^0_RM\to 0.
$$
Notice as well, that
\set\begin{equation}\label{eq_almost-vanish-LF}
(LF_R^jM)^a=0
\qquad
\text{in $\sD(R^a\Mod)$ for $j=1,\dots,p-1$}
\end{equation}
in view of \eqref{eq_decompose_F}, remark
\ref{rem_before-prooced}(iii), and proposition
\ref{prop_small-degrees}. The assertion now follows from
\eqref{eq_boundary-functor} and claim \ref{cl_systematic}.

(iii): According to \cite[I.4.3.1.7]{Il}, there is a natural
complex of functors :
$$
\Sigma
\quad:\quad
0\to G^d_R\xrightarrow{\ \partial_{d-1}}G_R^{d-1}\to\cdots\to
G_R^1\xrightarrow{\ \partial_0}G_R^0\to 0.
$$
Set $Z'_p:=\Ker\,\partial_{p-1}$; due to remark
\ref{rem_homogeneous_functor}(iii) we may consider the
derived functor $LZ'_p$, and since $\Sigma$ is exact on flat
$R$-modules, we obtain two exact sequences of $R$-modules :
$$
\begin{aligned}
\Sigma' \quad : \quad &
0\to LG^d_RM\to LG_R^{d-1}M\to\cdots\to LG^p_RM\to LZ'_pM\to 0 \\
\Sigma'' \quad : \quad &
0\to LZ'_pM\to LG_R^{p-1}M\to\cdots\to LG^1_RM\to LG^0_RM\to 0.
\end{aligned}
$$
On the one hand, in light of (ii), remark
\ref{rem_before-prooced}(iii,vi) and \eqref{eq_decompose_F},
we see that
$$
\cosk_{n+p-1}(LG^j_RM)^a=0
\qquad
\text{in $\sD(R^a\Mod)$ for every $j=1,\dots,d$}.
$$
By applying claim \ref{cl_systematic} to the exact sequence
$\cosk_{n+p-1}\Sigma'{}^a$, we deduce that $H_i(LZ'_pM)^a=0$
for every $i<n+p-1$. On the other hand, \eqref{eq_decompose_F}
and proposition \ref{prop_small-degrees} imply that
$$
(LG_R^jM)^a=0
\qquad
\text{in $\sD(R^a\Mod)$ for $j=1,\dots,p-1$}
$$
so the assertion follows, after applying claim \ref{cl_systematic}
to the exact sequence $\Sigma''$.
\end{proof}

\sset\subsubsection{}\label{subsec_attack}
Let now $\phi:R\to S$ be any morphism of simplicial
$V^a$-algebras; pick any resolution $\rho:P\to S$ with $P$
a flat $R$-algebra, and consider the composition
\set\begin{equation}\label{eq_will-be-Delta}
S\otimes_RP\xrightarrow{\ S\otimes_R\rho\ }
S\otimes_RS\xrightarrow{\ \mu_S\ }S
\end{equation}
where $\mu_S$ is the multiplication law of $S$. Then
\eqref{eq_will-be-Delta} represents a morphism
$$
\Delta(\phi):S\ellotimes_RS\to S
\qquad
\text{in $\sD(R\Alg)$}
$$
which is independent (up to unique isomorphism) of the
choice of $P$. Notice that if $\phi$ is an isomorphism
in $\sD(s.V^a\Alg)$, then the same holds for
$\phi\ellotimes_R\phi:R\to S\ellotimes_RS$; also, we have
$$
\phi=\Delta(\phi)\circ(\phi\ellotimes_R\phi)
\qquad
\text{in $\sD(s.V^a\Alg)$}.
$$
More generally, set $C:=\Cone\,\phi$; then, by considering
the section $S\ellotimes_R\phi$ of $\Delta(\phi)$, we get
a natural decomposition
$$
S\ellotimes_RS\isom S\oplus(S\ellotimes_RC)
\qquad
\text{in $\sD(R\Mod)$}
$$
which identifies $\Delta(\phi)$ with the natural projection,
whence a natural isomorphism
\set\begin{equation}\label{eq_heart}
\Cone\,\Delta(\phi)\isom\sigma S\ellotimes_RC.
\end{equation}
Thus, say that $C=\sigma^kC'$ in $\sD(R\Mod)$, for some
$R$-module $C'$ with $H_0C'\neq 0$; there follows a
distinguished triangle in $\sD(R\Mod)$
$$
\sigma^kC'\to\sigma^kS\ellotimes_RC'\to
\sigma^{2k}C'\ellotimes_RC'\to\sigma^{k+1}C'.
$$
Especially, if $k\geq 2$, then $H_kC=H_k(S\ellotimes_RC)$,
and we see that, in this case,
$\Cone\,\Delta(\phi)=\sigma^{k+1}C''$ in $\sD(R\Mod)$, for some
$C''$ such that $H_0C''\neq 0$. However, it is clear from
\eqref{eq_heart} that $\Cone\,\Delta^2(\phi)=\sigma^2C''$
for some object $C''$ of $\sD(R\Mod)$. We conclude that if
$\Delta^n(\phi)$ is an isomorphism in $\sD(R\Mod)$ for some
$n\geq 2$, then the same holds already for $\Delta^2(\phi)$.

\begin{definition} In the situation of \eqref{subsec_attack},
we say that $\phi$ is a {\em weakly \'etale morphism}, if
$\Delta^2(\phi)$ is an isomorphism in $\sD(R\Alg)$.
\end{definition}

\sset\subsubsection{}\label{subsec_change-of-base}
Let $\phi:R'\to R$ be any morphism of simplicial $V^a$-algebras,
and $S$, $T$ two $R$-algebras. Since the standard resolution
$F_R(S):=F_{R_{!!}}(S_{!!})^a$ of example \ref{ex_standard-for-alg}
is clearly functorial in both $R$ and $S$, we have a natural map
$$
F_{R'}(S)^\Delta\to F_R(S)^\Delta
$$
of simplicial $R'$-algebras, whence a natural morphism
\set\begin{equation}\label{eq_change-of-base}
S\ellotimes_{R'}T\to F_R(S)^\Delta\otimes_{R'}T\to S\ellotimes_RT.
\end{equation}

\begin{lemma}\label{lem_derived-wet}
In the situation of \eqref{subsec_change-of-base}, suppose
that $\phi$ is a quasi-isomorphism. Then the same holds for
\eqref{eq_change-of-base}.
\end{lemma}
\begin{proof} Let $\psi:F_{R'}(S)^\Delta\otimes_{R'}R\to F_R(S)^\Delta$
be the natural morphism; by construction, \eqref{eq_change-of-base}
factors as the composition of $\psi\otimes_R\one_T$ and the
natural isomorphism
$$
F_{R'}(S)^\Delta\otimes_{R'}T\isom
(F_{R'}(S)^\Delta\otimes_{R'}R)\otimes_RT.
$$
Since both $F_{R'}(S)^\Delta\otimes_{R'}R$ and $F_R(S)^\Delta$
are flat $R$-algebras, remark \ref{rem-prod-for-dr-a-alg}(ii)
then reduces to checking that $\psi$ is a quasi-isomorphism.
Since the natural maps $\beta:F_R(S)^\Delta\to S$ and
$\alpha:F_{R'}(S)^\Delta\to F_{R'}(\Delta)\otimes_{R'}R$ are
quasi-isomorphisms, it suffices to show that the composition
$\beta\circ\psi\circ\alpha$ is a quasi-isomorphism. But the
latter is none else than the standard resolution $F_{R'}(S)\to S$,
whence the lemma.
\end{proof}

\begin{proposition}\label{prop_wet-in-derived}
Let $\phi:R\to S$ and $\phi':R'\to S'$ be two morphisms of
simplicial $V$-algebras, and suppose we have a commutative
diagram in $\sD(R\Alg)$
\set\begin{equation}\label{eq_hot--commutes}
{\diagram R \ar[r]^-\phi \ar[d]_\psi & S \ar[d]^{\psi'} \\
           R' \ar[r]^-{\phi'} & S'
\enddiagram}
\end{equation}
where $\psi$ and $\psi'$ are isomorphisms (in $\sD(R\Alg)$).
Then $\phi^a$ is weakly \'etale if and only if the same holds
for $\phi'{}^a$.
\end{proposition}
\begin{proof} Denote by $\Hot(V\Alg)$ the homotopy category
of simplicial $V$-algebras, and recall that the multiplicative
system of quasi-isomorphisms in $\Hot(V\Alg)$ admits a right
calculus of fractions (\cite[I.3.1.8(ii)]{Il}). We begin with
the following special case :

\begin{claim}\label{cl_special-hot}
Suppose that $\psi$ and $\psi'$ are also morphisms of simplicial
$V$-algebras, and that \eqref{eq_hot--commutes} commutes in the
category $\Hot(V\Alg)$. Then the proposition holds.
\end{claim}
\begin{pfclaim} Indeed, in this situation it follows easily
from lemma \ref{lem_derived-wet} that $\phi^a$ (resp. $\phi'{}^a$)
is weakly \'etale if and only if the same holds for
$\psi'{}^a\circ\phi^a$ (resp. for $\phi'{}^a\circ\psi^a$),
so we are reduced to the case where $R=R'$, $S=S'$, and
$\phi,\phi':R\to S$ are two homotopic morphisms of simplicial
$V$-algebras. In this case, there exist morphisms of simplicial
$V$-algebras $\phi'':R\to S''$ and $d_0,d_1:S''\to S$, such that
$d_0$ and $d_1$ are quasi-isomorphisms, and $\phi=d_0\circ\phi''$,
$\phi'=d_1\circ\phi''$ (see \cite[I.2.3.2]{Il}). Then the claim
follows by applying repeatedly lemma \ref{lem_derived-wet}.
\end{pfclaim}

Next, there exist a simplicial $V$-algebra $S''$ and morphisms
$\beta:S''\to S$ and $\gamma:S''\to S'$ of simplicial
$V$-algebras, such that both $\beta$ and $\gamma$ are
quasi-isomorphisms, and $\psi'=\gamma\circ\beta^{-1}$ in
$\sD(s.V\Alg)$. Moreover, there exist morphisms
$\phi'':R''\to S''$ and $\psi'':R''\to R$ such that $\psi''$
is a quasi-isomorphism, and the resulting diagram
$$
\xymatrix{
R'' \ar[r]^-{\phi''} \ar[d]_{\psi''} & S'' \ar[d]^\beta \\ 
R \ar[r]^-\phi & S
}$$
commutes in $\Hot(V\Alg)$. By claim \ref{cl_special-hot},
we may then replace $\phi$ by $\phi''$ and $\psi'$ by
$\beta$, after which we may assume that $\psi'$ is a morphism
of simplicial $V$-algebras.

Likewise, we may find morphisms $\beta':R''\to R$ and
$\gamma':R''\to R'$ of simplicial $V$-algebras that are
quasi-isomorphisms, and such that $\psi=\gamma'\circ\beta'{}^{-1}$
in $\sD(s.V\Alg)$; again by lemma \ref{lem_derived-wet},
we see that $\phi^a$ is weakly \'etale if and only if the
same holds for $\phi^a\circ\beta^a$, so we may replace $\phi$
by $\phi\circ\beta$, $\psi$ by $\psi\circ\gamma$, and assume
from start that also $\psi$ is a map of simplicial $V$-algebras.

In this situation, by applying once again lemma
\ref{lem_derived-wet} we see that $\phi^a$ (resp. $\phi'{}^a$)
is weakly \'etale if and only if the same holds for
$\psi'{}^a\circ\phi^a$ (resp. for $\phi'{}^a\circ\psi^a$).
Therefore, we may assume from start that $R=R'$, $S=S'$, and
$\phi,\phi':R\to S$ are two maps of simplicial $V$-algebras
that represent the same morphism in $\sD(s.V\Alg)$. In
this case, we may find a morphism $\theta:R'''\to R$
of simplicial $V$-algebras, such that $\theta$ is a
quasi-isomorphism, and
$\phi\circ\theta,\phi'\circ\theta:R'''\to S$ are homotopic
maps; then the assertion follows from claim \ref{cl_special-hot}
(and again, from lemma \ref{lem_derived-wet}: details left
to the reader).
\end{proof}

\begin{corollary}\label{cor_wet-in-derived}
Let $\phi:R\to S$ and $\psi:S\to T$ be any two morphisms
of simplicial $V^a$-algebras. We have :
\begin{enumerate}
\item
If $\phi$ and $\psi$ are weakly \'etale, the same holds
for $\psi\circ\phi$.
\item
If $\phi$ and $\psi\circ\phi$ are weakly \'etale, the
same holds for $\psi$.
\item
If $\phi$ is weakly \'etale, and $R\to R'$ is any morphism
of simplicial $V^a$-algebras, then $\phi':=R'\ellotimes_R\phi$
is weakly \'etale.
\item
$\phi$ is weakly \'etale if and only if the same holds for
$\Delta(\phi)$.
\end{enumerate}
\end{corollary}
\begin{proof} (iv) is immediate from the discussion of
\eqref{subsec_attack}.

(iii): Endow $S':=R'\ellotimes_RS$ with the $R'$-algebra
structure deduced from its left tensor factor; then the
assertion is an immediate consequence of the following
more general

\begin{claim} There exists a commutative diagram in $\sD(R'\Alg)$
$$
\xymatrix{ S'\ellotimes_{R'}S' \ar[rr] \ar[rd]_{\Delta(\phi')} & &
R'\ellotimes_R(S\ellotimes_RS) \ar[ld]^{R'\ellotimes_R\Delta(\phi)} \\
& S'
}$$
whose horizontal arrow is an isomorphism.
\end{claim}
\begin{pfclaim} Pick any resolution $\rho:P\to S$ with $P$ a flat
$R$-algebra; then $R'\otimes_RP$ represents $R'\ellotimes_RS$,
and we have a natural isomorphism of $R'$-algebras
$$
(R'\otimes_RP)\otimes_{R'}(R'\otimes_RP)\isom R'\otimes_R(P\otimes_RP).
$$
In view of remark \ref{rem_loop-and-suspend}(i), the
source of this map represents $S'\ellotimes_{R'}S'$, and
the target represents $R'\ellotimes_R(S\ellotimes_RS)$.
Under these identifications, the morphism
$R'\ellotimes_R\Delta(\phi)$ becomes the map
$\one_{R'}\otimes_R(\mu_S\circ(\rho\otimes_R\rho))$ (notation
of \eqref{subsec_attack}), whereas $\Delta(\phi')$ is
represented by the multiplication map
$\mu_{R'\otimes_RP}=\one_{R'}\otimes_R\mu_P$ of $R'\otimes_RP$.
The claim follows straightforwardly.
\end{pfclaim}

(i): Pick a resolution $P\to S_{!!}$ with $P$ a flat $R_{!!}$-algebra,
and a resolution $Q\to T_{!!}$ with $Q$ a flat $P$-algebra; in
light of proposition \ref{prop_wet-in-derived} it suffices
to show that the resulting morphism $R\to Q^a$ is weakly
\'etale, so we may replace $S$ by $P^a$ and $T$ by $Q^a$, and
assume from start that both $\phi$ and $\psi$ are flat
morphisms. Due to (iv), it then suffices to check that
$\Delta(\psi\circ\phi):T\otimes_RT\to T$ is weakly \'etale;
however, the latter can be factored as the composition
\set\begin{equation}\label{eq_again-and-again}
T\otimes_RT\isom T\otimes_S(S\otimes_RS)\otimes_ST
\xrightarrow{\ T\otimes_S\Delta(\phi)\otimes_ST\ }
T\otimes_ST\xrightarrow{\ \Delta(\psi)\ } T
\end{equation}
and by (iv) the maps $\Delta(\psi)$ and $\Delta(\phi)$
are weakly \'etale, so the same holds for
$T\otimes_S\Delta(\phi)\otimes_ST$, in view of (iii).
We may therefore replace $\psi$ by $\Delta(\psi)$ and
$\phi$ by $T\otimes_S\Delta(\phi)\otimes_ST$, after
which, we may assume that $\Delta(\psi)$ is an isomorphism
in $\sD(s.A\Alg)$. In this case, the factorization
\eqref{eq_again-and-again} makes it clear that
$\Delta(\psi\circ\phi)$ is weakly \'etale, as stated.

(ii): Set $S':=S\ellotimes_RS$ and $T':=S\ellotimes_RT$;
we may factor $\psi$ as the composition
$$
S\xrightarrow{\ S\ellotimes_R(\psi\circ\phi)\ } T'
\xrightarrow{\ T'\ellotimes_{S'}\Delta(\phi)\ } T
$$
so the assertion follows from (i), (iii) and (iv).
\end{proof}

\begin{theorem}\label{th_vanishing-headache}
Let $\phi:R\to S$ be a morphism of simplicial $V$-algebras,
$c,p\in\N$ two integers, with $p$ a prime, and suppose that
$\phi^a$ is weakly \'etale. We have :
\begin{enumerate}
\item
If the ideal\/ $\fm\cdot H_0(S)$ of\/ $H_0(S)$ satisfies
condition $(\bB)$ of\/ \cite[\S2.1.6]{Ga-Ra}, then
$$
\L_{S/R}^a=0
\qquad
\text{in $\sD(S^a\Mod)$}.
$$
\item
If $S$ is a $\Z_{(p)}$-algebra and $H_i\L_{S/R}=0$ for every
$i<c$, then
$$
H_i\L_{S/R}^a=0
\qquad
\text{for every $i<c+2p-1$}.
$$
\end{enumerate}
\end{theorem}
\begin{proof} We begin with the following :

\begin{claim}\label{cl_wants-a-name}
We may assume that both $\Delta(\phi^a)$ and $H_0(\phi)$
are isomorphisms.
\end{claim}
\begin{pfclaim} To ease notation, set $S':=S\ellotimes_RS$,
and consider the sequence of morphisms
\set\begin{equation}\label{eq_ease-S-S}
S\xrightarrow{\ \one_S\ellotimes_R\phi\ }S'
\xrightarrow{\ \mu_B\ }S
\end{equation}
where $\mu_S$ is the composition of the multiplication
map $S\otimes_RS\to S$ and the natural map $S'\to S\otimes_RS$;
the transitivity triangle (\cite[III.2.1.2]{Il}) relative to
\eqref{eq_ease-S-S} yields a natural isomorphism
$$
\L_{S'/S}\isom\sigma S\otimes_{S'}\L_{S'/S}\isom
\sigma S\otimes_{S'}(S'\otimes_S\L_{S/R})
\qquad
\text{in $\sD(S\Mod)$}
$$
where the last isomorphism follows from the base change
theorem of \cite[III.2.2.1]{Il}. Thus, $\L_{S'/S}\simeq
\sigma\L_{S/R}$ in $\sD(S\Mod)$, so assertion (i) (resp.
(ii)) holds for $\phi$, provided it holds for $\Delta(\phi)$,
and then the claim follows from corollary \ref{cor_wet-in-derived}(iv).
\end{pfclaim}

Henceforth, we assume that both $\Delta(\phi^a)$ and
$H_0(\phi)$ are isomorphisms.
Let $P\to S$ be a resolution, with $P$ an $R$-algebra, such
that $P[n]$ is a free $R[n]$-algebra for every $n\in\N$; set
$P':=S\otimes_RP$, and recall that there are natural isomorphisms
$$
\L_{S/R}\isom S\otimes_P\Omega^1_{P/R}\isom
S\otimes_{P'}\Omega^1_{P'/S}
$$
where $\Omega^1_{P/R}$ denotes the flat simplicial $P$-module
such that $\Omega^1_{P/R}[n]:=\Omega_{P[n]/R[n]}$ (and likewise
for $\Omega^1_{P'/S}$), with faces and degeneracies deduced
from those of $P$ and $R$ (resp. those of $P'$ and $S$), in
the obvious way. In other words, set
$$
J:=\Ker(\mu_S:P'\to S)
$$
(with $\mu_S$ as in the proof of claim \ref{cl_wants-a-name});
then $\L_{S/R}\simeq J/J^2$ in $\sD(S\Mod)$, and notice that
\set\begin{equation}\label{eq_J_a-vanish}
J^a=0
\qquad
\text{in $\sD(S^a\Mod)$}
\end{equation}
since $\Delta(\phi^a)$ is an isomorphism.

\begin{claim}\label{cl_J-is-reg}
$J$ is a quasi-regular ideal, and $H_0J=0$.
\end{claim}
\begin{pfclaim} Since $H_0(\phi)$ is an isomorphism, it
is clear that $H_0J=0$. Next, for every $n\in\N$, the
$S[n]$-algebra $P'[n]$ is free, hence we reduce to showing
the following. Let $B$ be any ring, $C:=B[X_i~|~i\in I]$ any
free $B$-algebra (for any set $I$), and $f:C\to B$ any morphism
of $B$-algebras; then $\Ker\,f$ is a quasi-regular ideal of $C$.
However, for every $i\in I$, set $b_i:=f(X_i)$, and let
$g:C\isom C$ be the isomorphism of $B$-algebras such that
$g(X_i)=X_i-a_i$ for every $i\in I$; clearly $\Ker\,f$ is
a quasi-regular ideal if and only if the same holds for
$g^{-1}\Ker\,f$. But the latter is the ideal $(X_i~|~i\in I)$,
so the assertion follows from remark \ref{rem_augm-of-free}(iii).
\end{pfclaim}

(i): We need to show that $H_n(J/J^2)^a=0$ for every $n\in\N$,
and we shall argue by induction on $n$. The assertion for $n=0$
is clear from \eqref{eq_J_a-vanish}, in view of the exact
sequence
$$
H_nJ\to H_n(J/J^2)\to H_{n-1}J^2
\qquad
\text{for every $n\in\N$}.
$$
Hence, suppose that $n>0$, and the assertion is already known
for every degree $<n$. The same exact sequence reduces to
checking that $H_{n-1}(J^2)^a=0$. Now, from claim
\ref{cl_J-is-reg} and theorem \ref{th_Quillen}, we know that
$H_{n-1}J^n=0$. Therefore, by an easy induction, we
are further reduced to showing that $H_{n-1}(J^i/J^{i+1})^a=0$
for every $i=2,\dots,n-1$. However, on the one hand, proposition
\ref{prop_quasi-regular}(ii.a) says that the natural map
$$
L\Sym^i_S(J/J^2)\to J^i/J^{i+1}
$$
is an isomorphism in $\sD(S\Mod)$, for every $i\in\N$; on
the other hand, the inductive assumption and theorem
\ref{th_pluri-vanish}(i) imply that $H_j(L\Sym^iJ/J^2)^a=0$
for every $j\leq n$, whence the contention.

(ii): We show, by induction on $n$, that $H_n(J/J^2)^a=0$
for every $n<c+2p-1$. Arguing as in the foregoing, we may
assume that $c+2p-1>n>0$, and the sought vanishing is already
known in degrees $<n$, in which case we reduce to checking
that $H_{n-1}(L\Sym^i_SJ/J^2)^a=0$ for every $i=2,\dots,n-1$.
Set $M:=\cosk_n(s.\trunc_nJ/J^2)$; then $H_jM=0$ for every
$j<c$, and the inductive assumption implies that $M^a=0$ in
$\sD(S^a\Mod)$, so $H_{n-1}(L\Sym^i_SM)^a=0$, by theorem
\ref{th_prime-vanishing}(iii). Lastly, by corollary
\ref{cor_truncate-and-derive}(i), we see that
$$
H_{n-1}(L\Sym^i_SM)=H_{n-1}(L\Sym^i_SJ/J^2)
\qquad
\text{for every $i\in\N$}
$$
whence the contention.
\end{proof}

\begin{remark} Let $A\to B$ be any morphism of simplicial
rings such that the multiplication map
$\mu_B:B\otimes_AB\to B$ induces an isomorphism
$H_0(\mu_B):H_0(B)\otimes_{H_0(A)}H_0(B)\isom H_0(B)$.
Pick a resolution $\phi:P\to B$ such that $P[n]$ is
a free $A[n]$-algebra for every $n\in\N$, and set
$P':=B\otimes_AP$,
$J:=\Ker\,(\mu_B\circ(B\otimes_A\phi):P'\to B)$. The
assumption on $B$ implies that $H_0J=0$. On the other
hand, arguing as in the proof of claim \ref{cl_J-is-reg}
we see that $J$ is a quasi-regular ideal. Taking into
account proposition \ref{prop_quasi-regular}(ii) the
spectral sequence associated with the $J$-adic filtration
of $P'$ can be written as
$$
E^2_{pq}:=H_{p+q}(L\Sym^q_B(\L_{B/A}))\Rightarrow
\Tor^A_{p+q}(B,B).
$$
Though the $J$-adic filtration is not finite, using proposition
\ref{prop_convergence-simple}(i) and theorem \ref{th_Quillen}
one sees that this spectral sequence is convergent, and the
filtration on its abutment is finite. We will not use this result.
\end{remark}

\begin{corollary}\label{cor_vanishing-cot-wetale}
Let $\phi:A\to B$ be a morphism of\/ $V^a$-algebras,
$p\in\N$ a prime integer, and suppose that
\begin{enumerate}
\alphaenu
\item
$V$ is a $\Z_{(p)}$-algebra
\item
the induced morphism $s.\phi:s.A\to s.B$ of simplicial
$V^a$-algebras is weakly \'etale.
\end{enumerate}
Then $H_i\L_{B/A}^a=0$ for every $i\leq 2p$.
\end{corollary}
\begin{proof} Notice that condition (a) implies that
$B_{!!}$ is a $\Z_{(p)}$-algebra; letting $c:=0$ in
theorem \ref{th_vanishing-headache}(ii), we deduce
that $H_i\L_{B/A}^a=0$ for $i=0,1$. But in view of
claim \ref{cl_this-is-a}, the latter implies that
$H_i\L_{B_{!!}/A_{!!}}=0$ for $i=0,1$. So, actually the
morphism $\phi$ fulfills the conditions of theorem
\ref{th_vanishing-headache}(ii), with $c=2$, whence
the assertion.
\end{proof}

\begin{remark}\label{rem_instead-of-3.2.1}
As an application we can now generalize the deformation
theory for almost rings in \cite[\S3.2]{Ga-Ra} to the
case of an arbitrary basic setup $(V,\fm)$. Indeed,
\cite[Lemmata 3.2.1 and 3.26, Prop.3.2.9 and
Cor.3.2.11]{Ga-Ra} hold {\em verbatim} for any such
basic setup, with the same proof : the only difference
is that instead of invoking \cite[3.2.1]{Ga-Ra} we must
appeal to our theorem \ref{th_instead-of-3.2.1} in the
explanation of \cite[\S3.2.7]{Ga-Ra}. Next, let $A\to B$
be any weakly \'etale morphism of $V^a$-algebras; for every
prime integer $p$, let $V_{(p)}:=\Z_{(p)}\otimes_\Z V$, and
denote
$$
j_p:(V,\fm)^a\Alg\to(V_{(p)},\fm V_{(p)})^a\Alg
\qquad
C\mapsto\Z_{(p)}\otimes_\Z C
$$
the natural functor; with this notation, it is easily seen
that
$$
\Z_{(p)}\otimes_\Z B_{!!}=(j_pB)_{!!}
$$
(details left to the reader), and in view of \cite[II,2.2.1]{Il}
we deduce a natural isomorphism
$$
\Z_{(p)}\otimes_\Z\L_{B_{!!}/A_{!!}}\isom\L_{(j_pB)_{!!}/(j_pA)_{!!}}.
$$
Corollary \ref{cor_vanishing-cot-wetale} then easily implies
that $H_i\L^a_{B/A}=0$ for every $i\leq 4$. Combining with
the results just quoted from \cite{Ga-Ra}, we see that
the whole of \cite[Th.3.2.18]{Ga-Ra} still holds under
the current assumptions. In the same vein, we may complement
as follows the descent results of \cite[\S3.4]{Ga-Ra} :
\end{remark}

\begin{corollary}\label{cor_instead-of-3.2.1}
Consider a cartesian diagram of\/ $V^a$-algebras
$$
\cD
\qquad :\qquad
{\diagram A_0 \ar[r] \ar[d] & A_1 \ar[d] \\
A_2 \ar[r] & A_3
\enddiagram}$$
such that the morphism $A_1\to A_3$ is an epimorphism
on the underlying $V^a$-modules. Then the resulting
morphism $A_0\to A_1\times A_2$ is of universal effective
descent for the fibred categories $\mathbf{w.\acute{E}t}$
of weakly \'etale morphisms and $\mathbf{\acute{E}t}$ of
\'etale morphisms.
\end{corollary}
\begin{proof} Arguing as in the proof of
\cite[Prop.3.4.33]{Ga-Ra}, and using \cite[Cor.3.4.22]{Ga-Ra},
we see that the morphism $A_0\to A_1\times A_2$ is of effective
descent for the fibred categories $\mathbf{w.\acute{E}t}$ and
$\mathbf{\acute{E}t}$. Next, consider any morphism $A_0\to B_0$
of $V^a$-algebras, and set $B_i:=B_0\otimes_{A_0}A_i$ for
$i=1,2,3$. The resulting diagram $B_0\otimes_{A_0}\cD$ is
not necessarily cartesian; however, we notice :

\begin{claim}\label{cl_we-show-more-generally}
The induced morphism $f:B_0\to B'_0:=B_1\times_{B_3}B_2$
is an epimorphism on the underlying $V^a$-modules, with
nilpotent kernel.
\end{claim}
\begin{pfclaim} Indeed, we have a short exact sequence
of $A_0$-modules
$$
0\to A_0\to A_1\times A_2\to A_3\to 0
$$
inducing a right exact sequence
$B_0\to B_1\times B_2\to B_3\to 0$, which shows that
$f$ is an epimorphism. Newt, for any $B$-module $M$ we
show more generally that
$$
\Ann_B(M\otimes_{A_0}A_1)\cdot\Ann_B(M\otimes_{A_0}A_2)
\subset\Ann_B(M).
$$
Letting $M:=B$, it will then follow easily that $(\Ker\,f)^2=0$.
Now, let $x,y\in B_*$ be two almost elements such that
$x\cdot(M\otimes_{A_0}A_1)=0$ and $y\cdot(M\otimes_{A_0}A_2)=0$,
and let also $I:=\Ker\,(A_1\to A_3)=\Ker\,(A_0\to A_2)$;
since $M\otimes_{A_0}A_2=M/IM$, we deduce that
$yM\subset IM$, and on the other hand the natural morphism
of $A_0$-modules
$M\otimes_{A_0}A_1\otimes_{A_1}I\isom M\otimes_{A_0}I\to IM$
is an epimorphism, so that $x\cdot IM=0$, and the assertion
follows.
\end{pfclaim}

Now, by the foregoing we known already that the induced
morphism $B'_0\to B_1\times B_2$ is of effective
descent for $\mathbf{w.\acute{E}t}$ and $\mathbf{\acute{E}t}$.
But claim \ref{cl_we-show-more-generally} together with remark
\ref{rem_instead-of-3.2.1} implies that the base
change functors
$$
B_0\wEt\to B'_0\wEt
\qquad\text{and}\qquad
B_0\Et\to B'_0\Et
$$
are equivalences of categories; thus, also the induced morphism
$B_0\to B_1\times B_2$ is of effective descent, as required.
\end{proof}

\sset\subsubsection{}
We conclude this section with a few words on the
cohomology of sheaves of almost modules. Namely, let
$(V,\fm)$ again be an arbitrary basic setup, which we
view as a basic setup relative to the one-point topos
$\{\mathrm{pt}\}$, in the sense of \cite[\S3.3]{Ga-Ra}.
Let $(X,\cO_{\!X})$ be a ringed topos, $\fm_X\subset\cO_{\!X}$
an ideal, such that $(\cO_{\!X},\fm_X)$ is a basic setup
relative to the topos $X$, and suppose also that we are
given a morphism of ringed topoi :
$$
\pi:(X,\cO_{\!X})\to(\{\mathrm{pt}\},V)
\quad\text{such that}\quad
(\pi^*\fm)\cdot\cO_{\!X}\subset\fm_X.
$$
We deduce a natural morphism :
$$
R\pi_*:\sD^+((\cO_{\!X},\fm_X)^a\Mod)\to\sD^+((V,\fm)^a\Mod)
$$
which can be constructed as usual, by taking injective
resolutions. In case $\tilde\fm$ is a flat $V$-module, this
is the same as setting
$$
R\pi_*(K^\bullet):=(R\pi_*K^\bullet_!)^a.
$$
Using \cite[Cor.2.2.24]{Ga-Ra} one may verify that the two
definitions coincide : the details shall be left to the reader.
In many cases, both statements and proofs of results
concerning the cohomology of $\cO_{\!X}$-modules carry over
{\em verbatim} to $\cO^a_{\!X}$-modules. One sets, as customary :
$$
H^\bullet(X,K^\bullet):=H^\bullet R\pi_*K^\bullet
\qquad
\text{for every object $K^\bullet$ of $\sD^+((\cO_{\!X},\fm_X)^a\Mod)$}.
$$

\sset\subsubsection{}\label{subsec_descend-cohom}
As an illustration, we consider the following situation. Suppose
that :
\begin{itemize}
\item
$f:A\to A'$ is a map of $V$-algebras, and set :
$$
\phi:=\Spec\,f:X':=\Spec\,A'\to X:=\Spec\,A.
$$
\item
$t\in A$ is an element which is regular both in $A$ and
in $A'$, and such that the induced map $A/tA\to A'/tA'$
is an isomorphism.
\item
$U\subset X$ is any open subset containing
$D(t):=\{\fp\in X~|~t\notin\fp\}$, and set $U':=\phi^{-1}U$.
\item
$\cF$ is a quasi-coherent $\cO^a_{\!U}$-module, and
set $\cF':=\phi^*_{|U'}\cF$, which is a quasi-coherent
$\cO^a_{\!U'}$-module.
\end{itemize}
Then we have natural morphisms of $A^a$-modules :
\set\begin{equation}\label{eq_nothing-more}
H^q(U,\cF)\to H^q(U',\cF')\qquad\text{for every $q\in\N$.}
\end{equation}

\begin{lemma}\label{lem_compare-coh-t-reg}
In the situation of \eqref{subsec_descend-cohom}, suppose
moreover that $t$ is $\cF$-regular. Then \eqref{eq_nothing-more}
is an isomorphism for every $q>0$, and induces an isomorphism
of $A^{\prime a}$-modules :
$$
H^0(U,\cF)\otimes_AA'\isom H^0(U',\cF').
$$
\end{lemma}
\begin{proof} To ease notation, we shall write $\phi$ instead
of $\phi_{|U'}$. To start out, we remark :

\begin{claim}\label{cl_double-vanish}
(i)\ \ $\cTor_1^{\cO^a_{\!U}}(\cF,\cO^a_{\!U}/t\cO^a_{\!U})=0$
and $\Tor_1^{A^a}(H^0(U,\cF),A/tA)=0$.
\begin{enumerate}
\addenu
\item
The $\cO^a_{\!U}$-module
$\phi_*\cF'=\cF\otimes_{\cO^a_{\!U}}\phi_*\cO^a_{\!U'}$
(resp. the $A^{\prime a}$-module $H^0(U,\cF)\otimes_{A^a}A^{\prime a}$)
is $t$-torsion-free.
\end{enumerate}
\end{claim}
\begin{pfclaim}(i): Since $t$ is regular on $\cO_{\!U}$, we have
a short exact sequence :
$\cE:=(0\to\cO^a_{\!U}\to\cO^a_{\!U}\to\cO^a_{\!U}/t\cO^a_{\!U}\to 0)$,
and since $t$ is regular on $\cF$, the sequence
$\cE\otimes_{\cO^a_{\!U}}\cF$ is still exact; the vanishing of
$\cTor_1^{\cO^a_{\!U}}(\cF,\cO^a_{\!U}/t\cO^a_{\!U})$ is an
easy consequence. For the second stated vanishing one argues
similarly, using the exact sequence $0\to A\to A\to A/tA\to 0$.

(ii): Under the current assumptions, the natural map
$\cO_{\!U}/t\cO_{\!U}\to\phi_*(\cO_{\!U'}/t\cO_{\!U'})$ is an
isomorphism, whence a short exact sequence :
$\cE':=(0\to\phi_*\cO^a_{\!U'}\to\phi_*\cO^a_{\!U'}\to
\cO^a_{\!U}/t\cO^a_{\!U}\to 0)$. In view of (i), the sequence
$\cE'\otimes_{\cO^a_{\!U}}\cF$ is still exact, so
$\cF\otimes_{\cO^a_{\!U}}\phi_*\cO^a_{\!U'}$ is $t$-torsion-free.
An analogous argument works as well for
$H^0(U,\cF)\otimes_{A^a}A^{\prime a}$.
\end{pfclaim}

\begin{claim}\label{cl_thicky} For every $n>0$, the map
$A/t^nA\to A'/t^nA'$ induced by $f$ is an isomorphism.
\end{claim}
\begin{pfclaim} Indeed, since $t$ is regular on both $A$ and
$A'$, and $A/tA\isom A'/tA'$, the map of graded rings
$\oplus_{n\in\N}t^nA/t^{n+1}A\to\oplus_{n\in\N}t^nA'/t^{n+1}A'$
is bijective. Then the claim follows from
\cite[Ch.III, \S2, n.8, Cor.3]{BouAC}.
\end{pfclaim}

Let $j:D(t)\to U$ be the natural open immersion, and set :
$$
\cF[t^{-1}]:=j_*j^*\cF\simeq\cF\otimes_{\cO^a_{\!U}}j_*\cO^a_{\!D(t)}.
$$
Notice the natural isomorphisms :
\set\begin{equation}\label{eq_also-for-prime}
j_*j^*\phi_*\cF'\simeq
(\cF\otimes_{\cO^a_{\!U}}\phi_*\cO^a_{\!U'})
\otimes_{\cO^a_{\!U}}j_*\cO^a_{\!D(t)}
\simeq\cF[t^{-1}]\otimes_{\cO^a_{\!U}}\phi_*\cO^a_{\!U'}.
\end{equation}
Since $t$ is $\cF$-regular, the natural map $\cF\to\cF[t^{-1}]$
is a monomorphism, and the same holds for the corresponding
map $\phi_*\cF'\to j_*j^*\phi_*\cF'$, in view of claim
\ref{cl_double-vanish}(ii). Moreover, $\cG:=\cF[t^{-1}]/\cF$
can be written as the increasing union of its subsheaves
$\Ann_\cG(t^n)$ (for all $n\in\N$); hence claim \ref{cl_thicky}
implies that the natural map
$\cG\to\cG\otimes_{\cO^a_{\!U}}\phi_*\cO^a_{\!U'}$
is an isomorphism. There follows a ladder of short exact sequences :
\set\begin{equation}\label{eq_ladder-on-U}
{\diagram
0 \ar[r] & \cF \ar[r] \ar[d] & \cF[t^{-1}] \ar[r] \ar[d] &
\cG \ar[r] \ddouble & 0 \\
0 \ar[r] & \phi_*\cF' \ar[r] &
\cF[t^{-1}]\otimes_{\cO^a_{\!U}}\phi_*\cO^a_{\!U'} \ar[r] &
\cG \ar[r] & 0.
\enddiagram}
\end{equation}
On the other hand, since $j$ is an affine morphism, we may
compute :
$$
H^q(U,\cF[t^{-1}])\simeq H^q(U,Rj_*j^*\cF)\simeq H^q(D(t),j^*\cF)
\simeq 0 \qquad \text{for every $q>0$.}
$$
Likewise, from \eqref{eq_also-for-prime} we get :
$$
H^q(U,\cF[t^{-1}]\otimes_{\cO^a_{\!U}}\phi_*\cO^a_{\!U'})\simeq 0
\qquad \text{for every $q>0$.}
$$
Thus, in the commutative diagram :
$$
\xymatrix{
H^{q-1}(U,\cG) \ar[r]^-\partial \ar[d]_-{\partial'} &
H^q(U,\cF) \ar[d] \\
H^q(U,\phi_*\cF) \ar[r]^-\sim & H^q(U',\cF')
}$$
the boundary maps $\partial$ and $\partial'$ are isomorphisms
whenever $q>1$, and the right vertical arrow is \eqref{eq_nothing-more},
so the assertion follows already for every $q>1$.
To deal with the remaining cases with $q=0$ or $1$, we look at the
ladder of exact cohomology sequences:
$$
\xymatrix{
0 \ar[r] & H^0(U,\cF) \ar[r]^-\alpha \ar[d] &
H^0(D(t),\cF) \ar[r]^-\beta \ar[d] &
H^0(U,\cG) \ar[r] \ddouble & H^1(U,\cF) \ar[r] \ar[d] & 0 \\
0 \ar[r] & H^0(U',\cF') \ar[r] &
H^0(D(t),\cF)\otimes_{A^a}A^{\prime a} \ar[r]^-{\beta'} &
H^0(U,\cG) \ar[r] & H^1(U',\cF') \ar[r] & 0
}$$
deduced from \eqref{eq_ladder-on-U}.
On the one hand, we remark that the natural inclusion
$\Img\,\beta\subset\Img\,\beta'$ factors as a composition
$$
\Img\,\beta\xrightarrow{\iota}
M:=A^{\prime a}\otimes_{A^a}\Img\,\beta
\xrightarrow{\pi}\Img\,\beta'
$$
where $\iota$ is given by the rule : $x\mapsto 1\otimes x$,
for all $x\in A^a_*$, and $\pi$ is an epimorphism.

On the other hand, the image of $\beta$ is isomorphic to $N:=\Coker\,\alpha$,
and the latter is the increasing union of its submodules
$\Ann_N(t^n)$ (for all $n\in\N$). Again from claim \ref{cl_thicky}
we deduce that $M=\Img\,\beta$, hence $\Img\,\beta'=\Img\,\beta$,
which shows that \eqref{eq_nothing-more} is an isomorphism also
for $q=1$. By the same token,
$\Ker(\alpha\otimes_{A^a}\one_{A^{\prime a}})$ is the increasing
union of its $t^n$-torsion submodules (for all $n\in\N$), since
it is a quotient of $\Tor_1^{A^a}(\Img\,\beta,A^{\prime a})$;
hence $\Ker(\alpha\otimes_{A^a}\one_{A^{\prime a}})=0$, in view
of claim \ref{cl_double-vanish}(ii). This easily implies
the last assertion for $q=0$.
\end{proof}

\subsection{Inverse systems of almost modules}
The considerations of this section expand upon
\cite[\S2.3.14, \S2.4.1, Lemmata 2.3.15, 2.4.2 and 2.4.13]{Ga-Ra}.

\sset\subsubsection{}\label{subsec_limit-of-alm-struct}
Let $V$ be a ring, $(\fm_\lambda~|~\lambda\in\Lambda)$ a
filtered system of ideals of $V$ such that
$\fm_\lambda^2=\fm_\lambda$ for every $\lambda$, and set
$\fm:=\bigcup_{\lambda\in\Lambda}\fm_\lambda$. Then $(V,\fm)$
and $(V,\fm_\lambda)$ for every $\lambda\in\Lambda$ are
basic setups, and we have an obvious compatible system
of localization functors :
\set\begin{equation}\label{eq_mamma-mia}
{\diagram
(V,\fm)^a\Mod \ar[r]^-{\pi_\lambda} \ar[rd]_{\pi_\mu} &
(V,\fm_\lambda)^a\Mod \ar[d]^{\pi_{\lambda\mu}} \\
& (V,\fm_\mu)^a\Mod
\enddiagram}
\qquad
\text{for every $\lambda\geq\mu$ in $\Lambda$}.
\end{equation}

\begin{proposition}\label{prop_limit-of-alm-struct}
(i)\ \
The compatible system \eqref{eq_mamma-mia} induces
an equivalence of categories
$$
\omega:
(V,\fm)^a\Mod\isom\Pslim{\lambda\in\Lambda}(V,\fm_\lambda)^a\Mod.
$$

(ii)\ \
If\/ $\fm_\lambda$ fulfills condition $(\bB)$ of\/
\cite[\S2.1.6]{Ga-Ra} for every $\lambda\in\Lambda$, then
the same holds for $\fm$.
\end{proposition}
\begin{proof}(ii) is obvious. To show (i), we construct a
quasi-inverse $\tau$ for $\omega$ as follows. Recall that
an object of the above $2$-limit is a datum
$(M_\bullet,f_{\bullet\bullet}):=((M_\lambda~|~\lambda\in\Lambda),
(f_{\lambda\mu}:\pi_{\lambda\mu}(M_\lambda)\isom M_\mu~|~\lambda\geq\mu))$
consisting of $(V,\fm_\lambda)^a$-modules $M_\lambda$ for every
$\lambda\in\Lambda$, and isomorphisms $f_{\lambda\mu}$
of $(V,\fm_\mu)^a$-modules, for every $\lambda\geq\mu$ in
$\Lambda$. Such a datum is required moreover to satisfy
the compatibility condition :
$$
f_{\mu\nu}\circ a_{\mu\nu}(f_{\lambda\mu})=f_{\lambda\nu}
\qquad
\text{whenever $\lambda\geq\mu\geq\nu$}.
$$
However, the discussion of \cite[\S2.2.2]{Ga-Ra} allows to
describe such a datum more concretely as follows. Choose,
for every $\lambda\in\Lambda$ a $V$-module $M'_\lambda$ whose
image in $(V,\fm_\lambda)^a\Mod$ represents $M_\lambda$; the
isomorphism $f_{\lambda\mu}$ corresponds then to a unique
$V$-linear isomorphism
$$
f'_{\lambda\mu}:\tilde\fm_\mu\otimes_VM'_\lambda\isom
\tilde\fm_\mu\otimes_VM'_\mu
\qquad
\text{whenever $\lambda\geq\mu$}
$$
(with $\tilde\fm_\lambda:=\fm_\lambda\otimes_V\fm_\lambda$
for every $\lambda\in\Lambda$) such that
\set\begin{equation}\label{eq_translates-compat}
f'_{\mu\nu}\circ(\tilde\fm_\nu\otimes_Vf'_{\lambda\mu})=
\tilde\fm_\nu\otimes_Vf'_{\lambda\nu}
\qquad
\text{whenever $\lambda\geq\mu\geq\nu$}.
\end{equation}
Set $M_{\lambda!}:=\tilde\fm_\lambda\otimes_VM'_\lambda$; by
composing the inverse of $f'_{\lambda\mu}$ with the natural
map $\tilde\fm_\mu\otimes_VM'_\lambda\to M_{\lambda!}$, we
obtain a well defined $V$-linear map
$$
f''_{\lambda\mu}:M_{\mu!}\to M_{\lambda!}
\qquad
\text{whenever $\lambda\geq\mu$}
$$
and in light of \eqref{eq_translates-compat} we deduce
easily that $f''_{\lambda\nu}=f''_{\mu\nu}\circ f''_{\lambda\mu}$
whenever $\lambda\geq\mu\geq\nu$ (details left to the
reader). We now let
$$
\tau(M_\bullet,f_{\bullet\bullet}):=
\bigl(\colim_{\lambda\in\Lambda}M_{\lambda!}\bigr)^a
\in\Ob((V,\fm)^a\Mod).
$$
If $(N_\bullet,g_{\bullet\bullet})$ is another object of
the above $2$-limit of categories, a morphism
$(M_\bullet,f_{\bullet\bullet})\to(N_\bullet,g_{\bullet\bullet})$
is the datum of a system of morphisms
$h_\lambda:M_\lambda\to N_\lambda$ fulfilling obvious compatibility
conditions with respect ot the isomorphisms $f_{\bullet\bullet}$
and $g_{\bullet\bullet}$; after choosing representative $V$-modules
$N'_\lambda$ for each $N_\lambda$, the morphism $h_\lambda$
corresponds to a unique $V$-linear map
$h'_\lambda:M_!\to N_!:=\tilde\fm_\lambda\otimes_VN'_\lambda$,
for every $\lambda\in\Lambda$. It is then clear that
the resulting system $h'_\bullet$ yields a natural transformation
$M_{\bullet!}\to N_{\bullet!}$ of direct systems indexed by
$\Lambda$, whence a well defined morphism of $(V,\fm)^a$-modules
$$
\tau(h_\bullet):=\bigl(\colim_{\lambda\in\Lambda}h'_\lambda\bigr)^a:
\tau(M_\bullet,f_{\bullet\bullet})\to\tau(N_\bullet,g_{\bullet\bullet}).
$$
This completes the construction of our functor $\tau$.
It is easily seen that the definition of $\tau$ is independent
of all choices, up to a canonical isomorphism of functors.
Especially, if $(M_\bullet,f_{\bullet\bullet})=\omega(M^a)$ for
some $V$-module $M$, we may assume that $M'_\lambda$ has
been chosen equal to $M$, and
$f'_{\lambda\mu}=\tilde\fm_\lambda\otimes_V\one_M$ for every
$\lambda,\mu\in\Lambda$ with $\lambda\geq\mu$. With these
choices, we get a natural isomorphism
$\tau\circ\omega(M^a)\isom M^a$ for every $(V,\fm)^a$-module
$M^a$. Next, for any given datum $(M_\bullet,f_{\bullet\bullet})$
as in the foregoing, set $(P_\bullet,t_{\bullet\bullet}):=
\omega\circ\tau(M_\bullet,f_{\bullet\bullet})$; for every
$\mu\in\Lambda$, the subset $\Lambda(\mu):=
\{\lambda\in\Lambda~|~\lambda\geq\mu\}$ is cofinal in
$\Lambda$, hence $P_\mu$ is naturally isomorphic to
$$
\bigl(\colim_{\lambda\in\Lambda(\mu)}M_{\lambda!}\bigr)^a\isom
\colim_{\lambda\in\Lambda(\mu)}(M_{\lambda!})^a\isom
\colim_{\lambda\in\Lambda(\mu)}\pi_{\lambda\mu}(M_\lambda)
$$
where the transition morphisms in the latter colimit
are given by the isomorphisms $f_{\lambda\mu}$; {\em i.e.}
is naturally isomorphic to $M_\mu$. Lastly, under this
identification, it is easily seen that the morphism
$t_{\lambda\mu}$ corresponds to $f_{\lambda\mu}$ whenever
$\lambda\geq\mu$ in $\Lambda$ (details left to the
reader); summing up, we have obtained a natural isomorphism
$\omega\circ\tau(M_\bullet,f_{\bullet\bullet})\isom
(M_\bullet,f_{\bullet\bullet})$, and the proof is concluded.
\end{proof}

\begin{corollary}\label{cor_limit-of-alm-struct}
With the notation of \eqref{subsec_limit-of-alm-struct}, for
every $(V,\fm)^a$-module $M$ and every $\lambda\in\Lambda$,
denote by $M_\lambda$ the image of $M$ in $(V,\fm_\lambda)^a\Mod$.
Let $A$ be a $(V,\fm)^a$-algebra, $M$ an $A$-module, $B$ an
$A$-algebra, and $r\in\N$. Then the following holds :
\begin{enumerate}
\item
We have natural isomorphisms of\/ $V$-modules :
$$
\colim_{\lambda\in\Lambda}M_{\lambda!}\isom M_!
\qquad
M_*\isom\lim_{\lambda\in\Lambda}M_{\lambda*}.
$$
\item
The system of localization functors $A\Alg\to A_\lambda\Alg$
induces a natural equivalence
$$
A\Alg\isom\Pslim{\lambda\in\Lambda} A_\lambda\Alg.
$$
\item
The $A$-module $M$ is flat (resp. faithfully flat, resp. almost
projective, resp. almost finitely generated, resp. almost finitely
presented) if and only if the same holds for the $A_\lambda$-module
$M_\lambda$, for every $\lambda\in\Lambda$.
\item
Suppose that $\fm_\lambda$ fulfills condition $(\bB)$ for
every $\lambda\in\Lambda$. Then the $A$-module $M$ is almost
projective of almost finite rank (resp. of finite rank $\leq r$)
if and only if the same holds for the $A_\lambda$-module $M_\lambda$,
for every $\lambda\in\Lambda$.
\item
The $A$-algebra $B$ is flat (resp. faithfully flat, resp.
weakly unramified, resp. weakly \'etale, resp. unramified,
resp. \'etale) if and only if the same holds for the
$A_\lambda$-algebra $B_\lambda$, for every $\lambda\in\Lambda$.
\end{enumerate}
\end{corollary}
\begin{proof}(i): The assertion for $M_!$ is clear. For the
assertion concerning $M_*$, it suffices to notice that for
any two $(V,\fm)^a$-modules $M$ and $N$, proposition
\ref{prop_limit-of-alm-struct}(i) implies that the natural map
\set\begin{equation}\label{eq_forza}
\Hom_{(V,\fm)^a\Mod}(M,N)\to\lim_{\lambda\in\Lambda}
\Hom_{(V,\fm_\lambda)^a\Mod}(M_\lambda,N_\lambda)
\end{equation}
is an isomorphism of $V$-modules.

(ii),(iii): Clearly the functors $\pi_\lambda$ and
$\pi_{\lambda\mu}$ of \eqref{subsec_limit-of-alm-struct} are
all compatible with tensor products and with the $\Alhom$
functors; moreover, for every $\lambda\in\Lambda$, every
short exact sequence of $(V,\fm_\lambda)^a$-modules is isomorphic
to $\Sigma^a$, for some short exact sequence $\Sigma$ of
$V$-modules (and likewise for short exact sequences of
$(V,\fm)^a$-modules : details left to the reader). In view
of proposition \ref{prop_limit-of-alm-struct}(i), assertion
(ii) follows straightforwardly, and taking into account the
isomorphism \eqref{eq_forza}, we also deduce that the
$A$-module $M$ is flat (resp. almost projective) if and only
if the same holds for the $A_\lambda$-module $M_\lambda$, for
every $\lambda\in\Lambda$.

Next, if every $M_\lambda$ is a faithfully flat $A_\lambda$-module,
and $X$ is an $A$-module with $M\otimes_AX=0$, it follows that
$M_\lambda\otimes_{A_\lambda}X_\lambda=0$ for every $\lambda\in\Lambda$,
whence $X_\lambda=0$ for every such $\lambda$, and hence $X=0$,
by proposition \ref{prop_limit-of-alm-struct}(i). Conversely,
if $M$ is faithfully flat and $\lambda\in\Lambda$, consider
any $A$-module $X$ such that $M_\lambda\otimes_{A_\lambda}X_\lambda=0$;
then $(M\otimes_AX)_\lambda=0$, {\em i.e.}
$\tilde\fm_\lambda\otimes(M_*\otimes_{A_*}X_*)=0$, so that
$M\otimes_A(\tilde\fm_\lambda\otimes_VX)=0$. By assumption,
this implies that $\tilde\fm_\lambda\otimes_VX=0$, and therefore
$X_\lambda=0$; this shows that $M_\lambda$ is faithfully flat.

Next, say that $A=R^a$ and $M=N^a$ for some $V$-algebra $R$
and some $R$-module $N$; suppose that each $M_\lambda$ is
almost finitely generated (resp. almost finitely presented),
and let $\fm_0\subset\fm$ be any finitely generated subideal.
We then find $\lambda\in\Lambda$ such that
$\fm_0\subset\fm_\lambda$, and by applying
\cite[Cor.2.3.13]{Ga-Ra} to the $A_\lambda$-module $M_\lambda$
we find a finitely generated (resp. finitely presented)
$R$-module $N'$ and an $R$-linear map $N'\to N$ whose kernel
and cokernel are annihilated by $\fm_0$. Then, by applying
\cite[Cor.2.3.13]{Ga-Ra} to the $A$-module $M$, we deduce
the assertion.

(iv): If every $M_\lambda$ is almost projective of almost
finite rank, we already know that $M$ is almost projective,
and it remains only to check that $M$ is of almost finite
rank. To this aim, let $\eps\in\fm$ be any element, choose
$\lambda\in\Lambda$ such that $\eps\in\fm_\lambda^2$, and
write $\eps=\sum_{i=1}^n\eps_i\eps'_i$ for some $n\in\N$ and
$\eps_1,\eps'_1,\dots,\eps_n,\eps'_n\in\fm_\lambda$; if
$M_\lambda$ is of almost finite rank, we may then find
$j\in\N$ such that $\eps_i\cdot\Lambda^j_{A_\lambda}M_\lambda=0$
for $i=1,\dots,n$, in which case it follows easily that
$\eps_i\eps'_i\cdot\Lambda^j_AM=0$ for every $i=1,\dots,n$
(details left to the reader), whence the assertion. One
argues likewise in case each $M_\lambda$ has rank $\leq r$ :
the details shall be left to the reader.

(v) follows immediately from (i), (ii) and (iii).
\end{proof}

\begin{remark}\label{rem_prod_of-alm-structures}
(i)\ \
Let $V$ be a ring, $\fm,\fm'\subset V$ two ideals such that
$(V,\fm)$ and $(V,\fm')$ are both basic setups, and let
$\fm'':=\fm\cdot\fm'$; set as usual $\tilde\fm:=\fm\otimes_V\fm$,
and define likewise the $V$-modules $\tilde\fm'$ and $\tilde\fm''$.
It is easily seen that $(V,\fm'')$ is a basic setup as well
Moreover, we claim that there exists a natural isomorphism
of $V$-modules :
\set\begin{equation}\label{eq_JEFTA}
\tilde\fm\otimes_V\tilde\fm'\isom\tilde\fm''
\qquad
x\otimes y\otimes x'\otimes y'\mapsto xx'\otimes yy'.
\end{equation}
Indeed, let $\mu:\fm\otimes_V\fm'\to\fm''$ be the multiplication
map; we know already that $\Ker\,\mu$ is annihilated by both
$\fm$ and $\fm'$, hence
$\fm''\otimes_V\mu:\fm''\otimes_V\fm\otimes_V\fm'\to\tilde\fm''$
is an isomorphism, and likewise for
$\mu\otimes_V\fm\otimes_V\fm':\tilde\fm\otimes_V\tilde\fm'\to
\fm''\otimes_V\fm\otimes_V\fm'$, so the same holds for their
composition, which is \eqref{eq_JEFTA}.

(ii)\ \
Let $M$ be any $V$-module; denote by $(M,\fm)^a\in(V,\fm)^a\Mod$
the image of $M$, and define likewise $(M,\fm')^a\in(V,\fm')^a\Mod$
and $(M,\fm'')^a\in(V,\fm'')^a\Mod$. We deduce from (i) natural
$V$-linear identifications
$$
\begin{array}{c}
(M,\fm'')^a_*\isom((M,\fm)^a_*,\fm')^a_*\isom((M,\fm')^a_*,\fm)^a_*\ \\
(M,\fm'')^a_!\isom((M,\fm)^a_!,\fm')^a_!\isom((M,\fm')^a_!,\fm)^a_!.
\end{array}
$$
\end{remark}

\sset\subsubsection{}
Let $V$ be any ring; a {\em $V$-linear abelian category} is
the datum of an abelian category $\cA$ together with a
$V$-module structure on $\Hom_\cA(A,B)$ for every $A,B\in\Ob(\cA)$,
such that the composition law of $\cA$ is a $V$-bilinear map
$$
\Hom_\cA(A,B)\times\Hom_\cA(B,C)\to\Hom_\cA(A,C)
\qquad
\text{for every $A,B,C\in\Ob(\cA)$}.
$$
We say that a functor $F:\cA\to\cA'$ between $V$-linear
abelian categories is {\em $V$-linear} if the map
$$
\Hom_\cA(A,B)\to\Hom_{\cA'}(FA,FB)
\qquad
f\mapsto Ff
$$
is $V$-linear, for every $A,B\in\Ob(\cA)$.

\begin{remark}\label{rem_V-lin-abelian-cat}
Let $V$ be a ring, and $\cA,\cA'$ two $V$-linear abelian
categories.

(i)\ \
For every $x\in V$ and $A\in\Ob(\cA)$ we define the suboobject
$xA:=\Img\,x\cdot\one_A\subset A$. Moreover, if $J\subset V$
is an ideal generated by a finite system $x_1,\dots,x_n\in V$,
we let $JA:=\sum_{i=1}^nx_iA\subset A$. Let us check that the
subobject $JA$ does not depend on the choice of a finite
system of generators for $J$ : it suffices to show that for
every $a\in J$ we have $aA\subset JA$; to this aim, say that
$a=\sum_{i=1}^na_ix_i$ for certain $a_1,\dots,a_n\in V$, and
for every $i=1,\dots,n$ denote by
$A\xrightarrow{e_i}A^{\oplus n}\xrightarrow{p_i}A$ the inclusion
into the $i$-the direct summand and the projection onto the
$i$-th direct factor of $A^{\oplus n}$. Set
$\phi:=\sum_{i=1}^nx_i\cdot p_i$ and $\psi:=\sum_{i=1}^na_i\cdot e_i$;
it is easily seen that $\Img\,\phi=JA$ and the $V$-bilinearity
of the composition law implies that $\phi\circ\psi=a\cdot\one_A$.
Since $\Img\,\phi\circ\psi\subset\Img\,\phi$, the assertion
follows.

(ii)\ \
Recall that the opposite $\cA^o$ of the abelian category $\cA$
is abelian; moreover, $\cA^o$ inherits from $\cA$ an obvious
$V$-linear structure : the reader may spell out the details.
Likewise, if $F:\cA\to\cA'$ any $V$-linear functor, then
$F^o:\cA^o\to\cA'^o$ is again $V$-linear.

(iii)\ \
Obviously, for every $V$-algebra $A$, the category $A\Mod$
is naturally a $V$-linear abelian category. If $(V,\fm)$
is any basic setup, and $A$ any $(V,\fm)^a$-algebra, then
also $A\Mod$ is a $V$-linear abelian category.

(iv)\ \
Let $J\subset V$ be any ideal, and $A\in\Ob(\cA)$ any object.
We say that $JA=0$ if $xA=0$ for every $x\in J$; notice that
if $J$ is finitely generated, this notation agrees with that
of (i). We say that $\cA$ is {\em $J$-torsion-free} if there
are no non-zero objects $A$ of $\cA$ with $JA=0$. For instance,
if $(V,\fm)$ is any basic setup, and $A$ any $(V,\fm)^a$-algebra,
then $A\Mod$ is $\fm$-torsion-free. Clearly, if $\cA$ is
$J$-torsion-free, the same holds for $\cA^o$, for its natural
$V$-linear structure.

(v)\ \
Let $f:A\to B$ be any morphism of $\cA$, and $J,J'\subset V$
two ideals such that $J\cdot\Ker\,f=J'\cdot\Coker\,f=0$ (notation
of (iv)). Then for every $a\in J$ and $b\in J'$ there exists
a morphism
$$
g:B\to A
\qquad\text{such that}\qquad
f\circ g=ab\cdot\one_B
\quad\text{and}\quad
g\circ f=ab\cdot\one_A.
$$
Indeed, write $f=e\circ p$, where
$A\xrightarrow{p}\Img\,f\xrightarrow{e}B$ are the natural
morphisms. Since $J\cdot\Ker\,p=J\cdot\Ker\,f=0$, the
endomorphism $a\cdot\one_A$ factors through $p$ and a
morphism $g_1:\Img\,f\to A$. Likewise, since
$J'\cdot\Coker\,e=J'\cdot\Coker\,f=0$, the endomorphism
$b\cdot\one_B$ factors through $e$ and a morphism
$g_2:B\to\Img\,f$. Then $g:=g_1\circ g_2$ fulfills the
stated conditions.

(vi)\ \
In the situation of (v), let $F:\cA\to\cA'$ be a $V$-linear
functor; we get $Ff\circ Fg=ab\cdot\one_{FB}$ and
$Fg\circ Ff=ab\cdot\one_{FA}$. It follows easily that
$JJ'\cdot\Ker\,Ff=JJ'\cdot\Coker\,Ff=0$.
\end{remark}

\begin{definition}\label{def_almost-Mittag-Leff}
Let $(V,\fm)$ be a basic setup, $\cA$ a $V$-linear abelian
category, $I$ a small category, and consider a functor
$$
M_\bullet:I\to\cA
\qquad
i\mapsto M_i
\qquad
(j\xrightarrow{\phi}i)\mapsto(M_j\xrightarrow{M_\phi}M_i).
$$

(i)\ \
We say that $M_\bullet$ is {\em essentially zero},
if for every $i\in\Ob(I)$ there exists a morphism
$\phi:j\to i$ in $I$ with $M_\phi=0$.

(ii)\ \
We say that $M_\bullet$ is {\em almost essentially zero},
if for every $i\in\Ob(I)$ and every subideal $\fm_0\subset\fm$
of finite type there exists a morphism $\phi:j\to i$ in $I$
with $\fm_0\cdot M_\phi=0$. We denote
$$
\bFun(I,\cA)_\mathrm{a.ess.0}
$$
the full subcategory of $\bFun(I,\cA)$ whose objects are
the almost essentially zero functors.

(iii)\ \
We say that $M_\bullet$ is {\em almost essentially constant},
if there exist $L\in\Ob(\cA)$ and a cone
$\pi_\bullet:c_L\Rightarrow M_\bullet$ inducing almost essentially
zero functors
$$
\Ker\,\pi_\bullet:I\to A\Mod
\qquad\text{and}\qquad
\Coker\,\pi_\bullet:I\to A\Mod.
$$

(iv)\ \
We say that $M_\bullet$ has the
{\em almost Mittag-Leffler property} if for every $i\in\Ob(I)$
and every subideal $\fm_0\subset\fm$ of finite type there exists
a morphism $j\xrightarrow{\phi}i$ in $I$ with :
$$
\fm_0\cdot\Img\,M_\phi\subset\Img\,M_{\phi\circ\psi}
\qquad
\text{for every morphism $\psi:k\to j$ in $I$}.
$$

(v)\ \
We say that $M_\bullet$ is a {\em Cauchy functor}, if for every
subideal $\fm_0\subset\fm$ of finite type there exists
$i\in\Ob(I)$ with
$\fm_0\cdot\Ker\,M_\phi=\fm_0\cdot\Coker\,M_\phi=0$ for every
morphism $\phi:j\to i$ of $I$.

We say that $M_\bullet$ is {\em null}, if it
is both a Cauchy functor and almost essentially zero.

(vi)\ \
Let $\bP(M_\bullet)$ be either one of the conditions : ``$M_\bullet$
is essentially zero'', ``$M_\bullet$ is almost essentially zero'',
``$M_\bullet$ is almost essentially constant'', ``$M_\bullet$ has
the almost Mittag-Leffler property'', ``$M_\bullet$ is a Cauchy
functor'' or ``$M_\bullet$ is null''. We say that
{\em $M_\bullet$ has the dual $\bP$ property} if $\bP(M^o_\bullet)$
holds, for the opposite functor $M^o_\bullet:I^o\to\cA^o$ (see
remark \ref{rem_V-lin-abelian-cat}(ii)).
\end{definition}

\begin{remark}\label{rem_almost-Mittag-Leff}
(i)\ \
In the situation of definition \ref{def_almost-Mittag-Leff},
for every $M\in\Ob(\cA)$ and every subideal $\fm_0\subset\fm$
of finite type, set $M[\fm_0]:=\bigcap_{x\in\fm_0}\Ker\,x\cdot\one_M$.
Notice that $M_\bullet$ is almost essentially zero if and only if the
following holds. For every such $\fm_0$, the induced functor
$$
M_\bullet/M_\bullet[\fm_0]:I\to\cA
\qquad
i\mapsto M_i/M_i[\fm_0]
$$
is essentially zero.

(ii)\ \
Suppose that $I$ is cofiltered and $\cA$ is complete, and
let us set :
$$
M^\triangle_i:=\bigcap_{\phi:j\to i}\Img\,M_\phi
\qquad
\text{for every $i\in\Ob(I)$}.
$$
We claim that $M_\psi$ restricts to a morphism in $\cA$
$$
M^\triangle_\psi:M^\triangle_{i'}\to M^\triangle_i
\qquad
\text{for every morphism $\psi:i'\to i$ in $I$}.
$$
Indeed, it suffices to check that
$M_\psi(M^\triangle_{i'})\subset\Img\,M_\phi$ for every morphism
$\phi:j\to i$ in $I$. However, since $I$ is cofiltered, for
every such $\phi$ there exist $k\in\Ob(I)$ and morphisms
$\phi':k\to i'$ and $\psi':k\to j$ such that
$\phi\circ\psi'=\psi\circ\phi'$. It follows that
$M_\psi(M^\triangle_{i'})\subset M_\psi(\Img\,\phi')=M_\phi(\Img\,\psi')$
whence the contention. Next, let
$\tau_\bullet:c_L\Rightarrow M_\bullet$ be any cone; notice that
\set\begin{equation}\label{eq_the-foregoing}
\Img\,\tau_i\subset M^\triangle_i
\qquad
\text{for every $i\in\Ob(I)$}.
\end{equation}
On the other hand, a direct inspection of the definitions shows that
$$
(M_\bullet/M^\triangle_\bullet)^\triangle_i=0
\qquad
\text{for every $i\in\Ob(I)$}.
$$
Hence, \eqref{eq_the-foregoing} implies that every cone
$c_{L'}\Rightarrow M_\bullet/M^\triangle_\bullet$ is the zero morphism,
so we must have $\lim_IM_\bullet/M^\triangle_\bullet=0$, and since
$\lim_I$ is a left exact functor (because it is a right
adjoint : see \eqref{sec_Fubini}) we conclude that the
inclusion $M^\triangle_\bullet\to M_\bullet$ induces an isomorphism :
$$
\lim_IM^\triangle_\bullet\isom\lim_IM_\bullet.
$$

(iii)\ \
In the situation of (ii), let $A$ be any $(V,\fm)^a$-algebra, and
take $\cA:=A\Mod$ (see remark \ref{rem_V-lin-abelian-cat}(iii));
notice then that $M_\bullet$ has the almost Mittag-Leffler property
if and only if the following holds. For every $i\in\Ob(I)$, the
induced system $M_{\bullet/i}:=(\Img\,M_\phi~|~\phi:j\to i)$ of
submodules of $M_i$ is a Cauchy net of the uniform space $\cI_A(M_i)$
defined as in \cite[Def.2.3.1]{Ga-Ra}. In this case, by
\cite[Lemma 2.3.5]{Ga-Ra}, the Cauchy sequence $M_{\bullet/i}$
admits a unique limit in $\cI_A(M_i)$, and by direct
inspection we see that this limit is $M^\triangle_i$, since
for every subideal $\fm_0\subset\fm$ of finite type there
exists a morphism $\phi:j\to i$ in $I$ such that
$\fm_0\cdot\Img\,M_\phi\subset M^\triangle_i$.

(iv)\ \
Notice that $M_\bullet:I\to\cA$ is a Cauchy functor if and
only if the following holds. For every subideal $\fm_0\subset\fm$
of finite type there exists $i\in\Ob(I)$ such that for every
pair of morphisms $k\xrightarrow{\psi}j\xrightarrow{\phi}i$
of $I$, we have $\fm_0\cdot\Ker\,M_\psi=\fm_0\cdot\Coker\,M_\psi=0$.
Indeed, the latter condition clearly implies that $M_\bullet$
is a Cauchy functor; conversely, suppose that $M_\bullet$ is
a Cauchy functor, and let $\fm_0\subset\fm_1\subset\fm$ be
subideals of finite type with $\fm_0\subset\fm_1^2$; by
assumption, we may find $i\in\Ob(I)$ such that, for every
$\phi$ and $\psi$ as in the foregoing, $\fm_1$ annihilates
all but the first and the fourth terms in the induced exact
sequence :
$$
0\to\Ker\,M_\psi\to\Ker\,M_{\phi\circ\psi}\to\Ker\,M_\phi\to
\Coker\,M_\psi\to\Coker\,M_{\phi\circ\psi}\to\Coker\,M_\phi\to 0.
$$
It follows easily that
$\fm_1\cdot\Ker\,M_\psi=\fm_0\cdot\Coker\,M_\psi$, as required.
Notice also that every Cauchy functor has the almost
Mittag-Leffler property : the details shall be left to the reader.

(v)\ \
Suppose again that $I$ is cofiltered. In view of (iv), we
deduce that a functor $M_\bullet:I\to\cA$ is null if
and only if the following holds. For every subideal
$\fm_0\subset\fm$ of finite type there exists $i\in\Ob(I)$
such that $\fm_0M_j=0$ for every morphism $\phi:j\to i$ :
the details shall be left to the reader. It follows easily
that the full subcategory of $\bFun(I,\cA)$ whose objects
are the null functors, is a Serre subcategory : again,
we leave the details to the reader.
\end{remark}

\begin{lemma}\label{lem_a.ess.0-is-Serre}
In the situation of definition {\em\ref{def_almost-Mittag-Leff}},
the following holds :

{\em (i)}\ \
The category $\bFun(I,\cA)_\mathrm{a.ess.0}$ is a Serre
subcategory of\/ $\bFun(I,\cA)$. Hence
$$
\cL(I,\cA):=\bFun(I,\cA)/\bFun(I,\cA)_\mathrm{a.ess.0}.
$$
is an abelian category whose objects are the functors
$I\to\cA$.

{\em (ii)}\ \
If $\cA$ is $\fm$-torsion-free (see remark
{\em\ref{rem_V-lin-abelian-cat}(iv)}), then for every
almost essentially zero functor $M_\bullet:I\to A\Mod$,
we have $\lim_IM_\bullet=0$.

{\em (iii)}\ \
Suppose that $I$ is cofiltered, and $\cA$ is complete
and $\fm$-torsion-free. Then the functor
$$
\Lim_I:\bFun(I,\cA)\to\cA
$$
(see \eqref{sec_Fubini}) factors through the localization
$\bFun(I,\cA)\to\cL(I,\cA)$ and a functor
$$
\cL(I,\cA)\to\cA
\qquad
M_\bullet\mapsto\lim_IM_\bullet.
$$
\end{lemma}
\begin{proof}(i): Consider a short exact sequence of functors
from $I$ to $A\Mod$ :
\set\begin{equation}\label{eq_short-ex-seq-over-I}
0\to M'_\bullet\xrightarrow{f_\bullet}
M_\bullet\xrightarrow{g_\bullet}M''_\bullet\to 0.
\end{equation}
It is easily seen that if $M_\bullet$ is almost essentially
zero, the same holds for $M'_\bullet$ and $M''_\bullet$. Conversely,
suppose that $M'_\bullet$ and $M''_\bullet$ are almost essentially zero,
and let $i\in\Ob(I)$, and $\fm_0\subset\fm$ a finitely generated
subideal. By assumption there exists a morphism $\phi:j\to i$
in $I$ such that $\fm_0\cdot\Img\,M'_\phi=0$, and also a
morphism $\psi:k\to j$ such that $\fm_0\cdot\Img\,M''_\psi=0$.
It follows that there exists a morphism $h:\fm_0M_k\to M'_j$
such that $f_j\circ h:\fm_0M_k\to M_j$ is the restriction of
$M_\psi$. Then $M_{\phi\circ\psi}(\fm_0^2M_k)=
\fm_0\cdot\Img(M_\phi\circ f_j\circ h)=
\fm_0\cdot\Img(f_i\circ M'_\phi\circ h)=0$, which shows that
$M_\bullet$ is almost essentially zero.

(ii): Suppose first that $M_\bullet$ is essentially zero. Let
$X\in\Ob(\cA)$ be any object, and
$\tau_\bullet:c_X\Rightarrow M_\bullet$ any cone; by assumption,
for evey $i\in\Ob(I)$ there exists a morphism $\phi:j\to i$
in $I$ with $M_\phi=0$. Then $\tau_i=M_\phi\circ\tau_j=0$ for
every such $i$, so that $\tau_\bullet=0$, whence the contention.

For a general almost essentially zero functor $M_\bullet$,
and every subideal $\fm_0\subset\fm$ of finite type, we
have a short exact sequence of functors
$$
0\to M_\bullet[\fm_0]\to M_\bullet\to M_\bullet/M_\bullet[\fm_0]\to 0
$$
whose right-most term is essentially zero, by remark
\ref{rem_almost-Mittag-Leff}(i)), so that
$\lim_IM_\bullet/M_\bullet[\fm_0]=0$ by the previous case.
On the other hand, clearly $\fm_0$ annihilates
$\lim_IM_\bullet[\fm_0]$.
The functor $\Lim_I$ is left exact (since it admits a left
adjoint : see \eqref{sec_Fubini}); it then follows that
$\fm_0$ annihilates as well $\lim_IM_\bullet$, and since
$\fm_0$ is arbitrary and $\cA$ is $\fm$-torsion-free, the
claim follows.

(iii): As $\Lim_I$ is left exact, in the situation of
\eqref{eq_short-ex-seq-over-I} we see easily from (ii)
that if $M''_\bullet$ is almost essentially zero, then
$f_\bullet$ induces an isomorphism
$\lim_IM'_\bullet\isom\lim_IM_\bullet$. Thus, it remains only to
check that if $M'_\bullet$ is almost essentially zero, then
$g_\bullet$ induces an isomorphism
$g:\lim_IM_\bullet\isom\lim_IM''_\bullet$. To this aim, we suppose
first that $M_\bullet$ is the constant functor of value $M$, and
that for every subideal $\fm_0\subset\fm$ of finite type there
exists $i\in\Ob(I)$ such that $\fm_0\cdot M'_i=0$.
Then, for every morphism $\phi:j\to i$ in $I$ we have as
well $\fm_0\cdot M'_j=0$. Hence, for such $\fm_0$ and $i$,
and for every $b\in\fm_0$ and every $\phi:j\to i$, the
endomorphism $b\cdot\one_M:M\to M$ is the composition of
$g_j:M\to M''_j$ and a unique morphism $h_j:M''_j\to M$ of
$\cA$, and it is easily seen that
$g_j\circ h_j=b\cdot\one_{M''_j}$. Moreover, the system of
morphisms $(g_j~|~\phi:j\to i)$ amounts to a natural cocone
$h_\bullet:M''_\bullet\circ\ss_i\Rightarrow M_\bullet\circ\ss_i$,
where $\ss_i:I/i\to I$ denotes the source functor. Since
$\ss_i$ is final (example \ref{ex_filtered-final}(i)),
we deduce a morphism of $\cA$
$$
h:L'':=\lim_IM''_\bullet\to M=\lim_IM_\bullet.
$$
By construction, we have $g\circ h=b\cdot\one_{L''}$ and
$h\circ g=b\cdot\one_M$. Especially, $b$ annihilates
$\Ker\,g$ and $\Coker\,g$, and since $\fm_0$ and $b$
are arbitrary, the assertion follows in this case.

Lastly, let $M'_\bullet$ be an arbitrary almost essentially
zero functor, and notice that
$$
N_\phi:=f_i(\Img\,M'_\phi)=M_\phi(\Ker\,g_j)
\qquad
\text{for every morphism $\phi:j\to i$ in $I$}.
$$
Hence, for every such $\phi$, the morphism $M''_\phi$ is
the composition of unique morphisms of $\cA$
$$
h_\phi:M''_j\xrightarrow{h_\phi}Q_\phi:=M_i/N_\phi
\xrightarrow{\bar g_\phi}M''_i
\qquad\text{such that}\qquad
h_\phi\circ g_j=\pi_\phi\circ M_\phi
\quad\text{and}\quad
\bar g_\phi\circ\pi_\phi=g_i
$$
where $\pi_\phi:M_i\to Q_i$ is the projection.
For any fixed $i\in\Ob(I)$, let $c_{M_i}:I/i\to\cA$ be the
constant functor of value $M_i$; we get a short exact sequence
of functors
$$
0\to N_\bullet\to c_{M_i}\xrightarrow{\pi_\bullet}Q_\bullet\to 0 
$$
where $N_\bullet,Q_\bullet:I/i\to\cA$ are given by the rules :
$\phi\mapsto N_\phi$ and respectively $\phi\mapsto Q_\phi$ for
every morphism $\phi:j\to i$ in $I$. By construction, for
every subideal $\fm_0\subset\fm$ of finite type there exists
$\phi\in\Ob(I/i)$ such that $\fm_0\cdot N_\phi=0$; by the
foregoing case, it follows that $\pi_\bullet$ induces an
isomorphism
$$
\pi_i:M_i\isom Q:=\lim_{I/i}Q_\bullet.
$$
On the other hand, the system $(h_\phi~|~\phi\in\Ob(I/i))$
amounts to a natural transformation
$h_\bullet:M''_\bullet\circ\ss_i\Rightarrow Q_\bullet$ (the details
are left to the reader), whence an induced morphism of $\cA$
$$
h_i:L'':=\lim_IM''_\bullet\isom\lim_{I/i}M''_\bullet\circ\ss_i\to Q.
$$
In turn, the system $(\pi_i^{-1}\circ h_i~|~i\in\Ob(I))$ amounts
to a cone $c_{L''}\Rightarrow M_\bullet$, whence an induced morphism
of $\cA$
$$
h:L''\to L:=\lim_IM_\bullet.
$$

\begin{claim}\label{cl_get-splitting-for-g}
$g\circ h=\one_{L''}$.
\end{claim}
\begin{pfclaim} Indeed, let also
$\tau_\bullet:c_L\Rightarrow M_\bullet$ and
$\tau''_\bullet:c_{L''}\Rightarrow M''_\bullet$ be the universal
cones; it suffices to show that $\tau''_i\circ g\circ h=\tau''_i$
for every $i\in\Ob(I)$. But we have
$$
\tau''_i\circ g\circ h=\tau''_i=g_i\circ\tau_i\circ h=
g_i\circ\pi_i^{-1}\circ h_i=
\bar g_\phi\circ\pi_\phi\circ\pi_i^{-1}\circ h_i
$$
for every morphism $\phi:j\to i$ in $I$, and on the other hand
$\tau''_i=M''_\phi\circ\tau''_j=\bar g_\phi\circ h_\phi\circ\tau''_j$
so we are reduced to checking that
$$
\pi_\phi\circ\pi_i^{-1}\circ h_i=h_\phi\circ\tau''_j
\qquad
\text{for every morphism $\phi:j\to i$ in $I$}.
$$
Now, let $\tau^Q_\bullet:c_Q\Rightarrow Q_\bullet$ be the universal
cone; then $\pi_\phi\circ\pi_i^{-1}=\tau^Q_\phi$ for every such
$\phi$, and finally, $\tau^Q_\phi\circ h_i=h_\phi\circ\tau''_j$,
as required.
\end{pfclaim}

From claim \ref{cl_get-splitting-for-g} we see that
$g$ is an epimorphism; but it is also a monomorphism, since
$\lim_IM'_\bullet=0$, and since the functor $\lim_I$ is left
exact. Hence $g$ is an isomorphism, as stated.
\end{proof}

\begin{remark}\label{rem_Cauchy-functors}
Keep the situation of definition \ref{def_almost-Mittag-Leff},
and suppose that $I$ is cofiltered.

(i)\ \
Consider a short exact sequence \eqref{eq_short-ex-seq-over-I}
of functors $I\to\cA$, such that $M'_\bullet$ and $M''_\bullet$
are Cauchy functors. Then $M_\bullet$ is also a Cauchy functor.
Indeed, let $\fm_0\subset\fm_1\subset\fm$ be two subideals
of finite type with $\fm_0\subset\fm_1^2$. Since $I$ is cofiltered,
taking into account remark \ref{rem_almost-Mittag-Leff}(iv), we
find $i\in\Ob(I)$ such that $\fm_1$ annihilates the kernel and
cokernel of both $M'_\phi$ and $M''_\phi$, for every morphism
$\phi:j\to i$; by the snake lemma, it follows easily that
$\fm_0\cdot\Ker\,M_\phi=\fm_0\cdot\Coker\,M_\phi=0$ for every
such $\phi$, whence the contention.

(ii)\ \
Moreover, if $X_\bullet,Y_\bullet:I\to\cA$ are Cauchy functors,
and $f_\bullet:X_\bullet\to Y_\bullet$ is any morphism of functors,
then $\Ker\,f_\bullet$, $\Coker\,f_\bullet$ and $\Img\,f_\bullet$
are Cauchy functors. Indeed, let us check first the assertion
for $Z_\bullet:=\Img\,f_\bullet$ : since the latter is a subfunctor
of $Y_\bullet$, for every subideal $\fm_0\subset\fm$ of finite
type we find $i\in\Ob(I)$ such that for every pair of morphisms
$k\xrightarrow{\psi}j\xrightarrow{\phi}i$ we have
$\fm_0\cdot\Ker\,Z_\psi=0$ (remark \ref{rem_almost-Mittag-Leff}(iv)).
Likewise, since $Z_\bullet$ is a quotient of $X_\bullet$, for every
such $\fm_0$ there exists $i'\in\Ob(I)$ such that for every pair
of morphisms $k'\xrightarrow{\psi'}j'\xrightarrow{\phi'}i'$ we
have $\fm_0\cdot\Coker\,Z_{\psi'}=0$. Since $I$ is cofiltered, we
find $i''\in\Ob(I)$ with morphisms $i\leftarrow i''\to i'$, and
then $\fm_0$ annihilates the kernel of cokernel of $Z_\lambda$,
for every morphism $\lambda:j''\to i''$ of $I$.

Lastly, by considering the short exact sequences of functors :
$$
0\to\Ker\,f_\bullet\to X_\bullet\to Z_\bullet\to 0
\qquad
0\to Z_\bullet\to Y_\bullet\to\Coker\,f_\bullet\to 0
$$
and combining with (i), we deduce the assertions for
$\Ker\,f_\bullet$ and $\Coker\,f_\bullet$.
\end{remark}

\begin{proposition}\label{prop_a-Mittag-Leffler}
With the notation of definition {\em\ref{def_almost-Mittag-Leff}}
and of remark {\em\ref{rem_almost-Mittag-Leff}(ii)}, assume that
$I$ is cofiltered, and $\cA$ is complete and $\fm$-torsion-free.
Then we have :

{\em (i)}\ \
If $M_\bullet:I\to\cA$ has the almost Mittag-Leffler
property, the following holds :
\begin{enumerate}
\alphaenu
\item
$M^\triangle_\phi:M^\triangle_j\!\to\!M^\triangle_i$ is an epimorphism
in $\cA$, for every morphism $\phi:j\!\to\!i$ of $I$.
\item
The inclusions $(M^\triangle_i\!\to\!M_i~|~i\in\Ob(I))$ define an
isomorphism $M^\triangle_\bullet\!\isom\!M_\bullet$ in $\cL(I,\cA)$.
\end{enumerate}

{\em (ii)}\ \
For every functor $M_\bullet:I\to\cA$ the following conditions
are equivalent :
\begin{enumerate}
\alphaenu
\item
$M_\bullet$ is almost essentially constant.
\item
$M_\bullet$ is isomorphic to a constant functor, in the category
$\cL(I,\cA)$ of lemma {\em\ref{lem_a.ess.0-is-Serre}(i)}.
\item
Every universal cone $c_L\Rightarrow M_\bullet$ is an isomorphism
in the category $\cL(I,\cA)$.
\item
$M_\bullet$ has the almost Mittag-Leffler property, and
$M^\triangle_\bullet$ is a Cauchy functor.
\end{enumerate}

{\em (iii)}\ \
Every Cauchy functor $M_\bullet:I\to\cA$ is almost essentially
constant.
\end{proposition}
\begin{proof}(i.a): Pick a morphism $k\xrightarrow{\psi}j$
such that $\fm_0\cdot\Img\,M_\psi\subset M^\triangle_j$ (remark
\ref{rem_almost-Mittag-Leff}(iii)). Then
$$
\fm_0M^\triangle_i\subset\fm_0\cdot\Img\,M_{\phi\circ\psi}=
M_\phi(\fm_0\cdot\Img\,M_\psi)\subset M^\triangle_\phi(M^\triangle_j).
$$
Since $\fm_0$ is arbitrary, the assertion follows.

(i.b): For every $i\in\Ob(I)$, let $\beta_i:M^\triangle_i\to M_i$
be the inclusion; we need to check that the resulting functor
$\Coker\,\beta_\bullet:I\to\cA$ is almost essentially zero.
Indeed, let $i\in\Ob(I)$ and $\fm_0\subset\fm$ any subideal
of finite type; since $M_\bullet$ has the almost Mittag-Leffler
property, there exists a morphism $\phi:j\to i$ in $I$ such
that $\fm_0\cdot\Img\,M_\phi\subset M^\triangle_i$. Hence, $\fm_0$
annihilates the image of the induced morphism
$\Coker\,\beta_j\to\Coker\,\beta_i$, whence the assertion.

(ii.a)$\Rightarrow$(ii.d): Let $L\in\Ob(\cA)$ with a
cone $\tau_\bullet:c_L\Rightarrow M_\bullet$ such that the
functors $\Ker\,\tau_\bullet$ and $\Coker\,\tau_\bullet$ are
almost essentially zero. This means that for every
$i\in\Ob(I)$ and every subideal $\fm_0\subset\fm$ of finite
type, there exists a morphism $\phi:j\to i$ such that
$\fm_0\cdot\Img\,M_\phi\subset\Img\,\tau_i$. But clearly
$\Img\,\tau_i\subset M^\triangle_i$, so we see that $M_\bullet$
has the almost Mittag-Leffler property. By the same token,
since $\fm_0M^\triangle_i\subset\fm_0\cdot\Img\,M_\phi$, we
get $\fm_0M^\triangle_i\subset\Img\,\tau_i$, and since $\fm_0$
is arbitrary, we conclude that $M^\triangle_i=\Img\,\tau_i$ for
every $i\in\Ob(I)$. Next, by assumption, for every $i\in\Ob(I)$
and every subideal $\fm_0\subset\fm$ of finite type there exists
a morphism $\phi:j\to i$ of $I$ such that $\fm_0$ annihilates
the image of the induced morphism $\Ker\,\tau_j\to\Ker\,\tau_i$;
but obviously the latter is a monomorphism, so the condition
means that for every such $\fm_0$ there exists $j\in\Ob(I)$
such that $\fm_0\cdot\Ker\,\tau_j=0$. Now, let $\psi:k\to j$ be
any morphism of $I$, and set $K:=\tau_k^{-1}(\Ker\,M^\triangle_\psi)$;
it is easily seen that  $K\subset\Ker\,\tau_j$, hence $\fm_0K=0$,
and finally $\fm_0\cdot\Ker\,M^\triangle_\psi=0$, since we have
already established that $\tau_k:L\to M^\triangle_k$ is an epimorphism.

(ii.d)$\Rightarrow$(ii.b):  In light of (i.b), we are reduced
to checking that the functor $M^\triangle_\bullet$ is almost essentially
constant. Thus, let $M^\triangle$ be the limit of the functor
$M^\triangle_\bullet$, and
$\tau_\bullet:c_{M^\triangle}\Rightarrow M^\triangle_\bullet$ the universal
cone; let also $\fm_0\subset\fm$ be any subideal of
finite type, and pick $i\in\Ob(I)$ such that
$\fm_0\cdot\Ker\,M^\triangle_\phi=0$ for every morphism $\phi:j\to i$
of $I$. Recall that the source functor $\ss_i:I/i\to I$ is cofinal,
so that $M^\triangle$ represents also the limit of
$M^\triangle_\bullet\circ\ss_i:I/i\to\cA$. It follows easily that
$\Ker(\tau_i:M^\triangle\to M^\triangle_i)$ represents the limit
of the induced functor
$$
(0\times_{M^\triangle_i}M^\triangle_\bullet)\circ\ss_i:I\to\cA
\qquad
(j\xrightarrow{\phi}i)\mapsto
\Ker\,M^\triangle_\phi=0\times_{M^\triangle_i}M^\triangle_j.
$$
But then it is clear that $\fm_0\cdot\Ker\,\tau_i=0$.
Especially, the functor $\Ker\,\tau_\bullet:I\to\cA$
is trivially almost essentially zero. To conclude, we show
that $\tau_i$ is an epimorphism for every $i\in\Ob(I)$.
To this aim, let $\fm_0\subset\fm$ be any subideal of finite
type; in view of (i.a), we may assume that
$\fm_0\cdot\Ker\,M^\triangle_\phi=0$ for every morphism $\phi:j\to i$.
Hence, for every $b\in\fm_0$ we have a unique morphism in $\cA$
$$
f_{\phi,b}:M^\triangle_i\to M^\triangle_j
\qquad\text{such that}\qquad
M^\triangle_\phi\circ f_{\phi,b}=b\cdot\one_{M^\triangle_i}.
$$
By inspecting the construction, it is easily seen that
$M^\triangle_\psi\circ f_{\phi\circ\psi,b}=f_{\phi,b}$ for every
morphism $\psi:k\to j$ in $I$, hence the rule
$\phi\mapsto f_{\phi,b}$ defines a cone
$c_{M^\triangle_i}\Rightarrow M^\triangle_\bullet\circ\ss_i$,
whence a unique morphism of $A$-modules
$$
h:M^\triangle_i\to M^\triangle
\qquad
\text{such that $f_{\phi,b}=\tau_j\circ h$ for every $\phi:j\to i$}.
$$
Especially, $\tau_i\circ h=M^\triangle_\phi\circ\tau_j\circ h=
M^\triangle_\phi\circ f_{\phi,b}=b\cdot\one_{M^\triangle_i}$,
and since $b\in\fm_0$ is arbitrary, the assertion follows.

(ii.b)$\Rightarrow$(ii.a): It suffices to prove :

\begin{claim} Let $M_\bullet,M'_\bullet:I\to\cA$ be two functors,
$f_\bullet:M_\bullet\Rightarrow M'_\bullet$ a natural transformation,
and suppose that $f_\bullet$ is an isomorphism in $\cL(I,\cA)$. Then
$M_\bullet$ is almost essentially constant if and only if the same
holds for $M'_\bullet$.
\end{claim}
\begin{pfclaim} Let $L$ and $L'$ be the limits of $M_\bullet$
and respectively $M'_\bullet$; denote by
$\tau_\bullet:c_L\Rightarrow M_\bullet$ and
$\tau'_\bullet:c_{L'}\Rightarrow M'_\bullet$ the respective universal
cones. There exists a unique morphism $f:L\to L'$ in $\cA$
that makes commute the diagram :
$$
\xymatrix{ c_L \ar[r]^-{c_f} \ar[d]_{\tau_\bullet} &
c_{L'} \ar[d]^{\tau'_\bullet} \\
M_\bullet \ar[r]^-{f_\bullet} & M'_\bullet.
}$$
But lemma \ref{lem_a.ess.0-is-Serre}(iii) implies that $f$ is
an isomorphism; as $f_\bullet$ is an isomorphism in $\cL(I,\cA)$,
it follows that $\tau_\bullet$ is an isomorphism in $\cL(I,\cA)$ if
and only if the same holds for $\tau'_\bullet$.
\end{pfclaim}

Lastly, clearly (c)$\Rightarrow$(a). Suppose then that (a)
holds, and let $\tau_\bullet:c_L\Rightarrow M_\bullet$ be a
universal cone; let also $\mu_\bullet:c_X\Rightarrow M_\bullet$ be
another cone that is an isomorphism in $\cL(I,\cA)$. Then
there exists a unique morphism $f:X\to L$ in $\cA$ such
that $\tau_\bullet\odot c_f=\mu_\bullet$. Set $\mu:=\lim_I\mu_\bullet$
and $\tau:=\lim_I\tau_\bullet$, so that $\tau\circ f=\mu$.
Now, $\tau$ is an isomorphism, and the same holds for
$\mu_\bullet$, by lemma \ref{lem_a.ess.0-is-Serre}(iii);
hence $f$ is an isomorphism, so the cone $\mu_\bullet$ is
universal, whence (c).

(iii): We know already that $M_\bullet$ has the almost Mittag-Leffler
property (remark \ref{rem_almost-Mittag-Leff}(iv)); moreover,
$\Ker\,M^\triangle_\phi\subset\Ker\,M_\phi$ for every morphism
$\phi$ of $I$, hence $M_\bullet$ fulfills condition (ii.c),
whence the contention.
\end{proof}

\begin{proposition}\label{prop_apply-F-V-linear}
Let $\cA$ and $\cA'$ be $V$-linear abelian categories, $I$
a small category, $F:\cA\to\cA'$ a $V$-linear functor, and
$M_\bullet,M'_\bullet:I\to\cA$ two functors. The following holds :

{\em(i)}\ \
Let $f_\bullet:M_\bullet\Rightarrow M'_\bullet$ be a natural
transformation. If $\Ker\,f_\bullet$ and $\Coker\,f_\bullet$
are almost essentially zero (resp. null, resp.
dual almost essentially zero, resp. dual null),
then the same holds for the kernel and cokernel of
$Ff_\bullet:FM_\bullet\Rightarrow FM'_\bullet$.

{\em(ii)}\ \
If $M_\bullet$ is a Cauchy functor, the same holds for
$FM_\bullet$.

{\em(iii)}\ \
If $M_\bullet$ is almost essentially constant (resp. dual
almost essentially constant), so is $FM_\bullet$.

{\em(iv)}\ \
Suppose that $F$ is right exact (resp. left exact). If $M_\bullet$
has the almost Mittag-Leffler property (resp. the dual almost
Mittag-Leffler property), then the same holds for $FM_\bullet$.
\end{proposition}
\begin{proof}(ii): The assertion follows immediately from
remark \ref{rem_V-lin-abelian-cat}(vi).

(i): We can write $f_\bullet=e_\bullet\circ p_\bullet$,
where $e_\bullet:I\to\cA$ is a monomorphism of functors, and
$p_\bullet:I\to\cA$ is an epimorphism of functors, and then
clearly it suffices to check the stated assertions for $e_\bullet$
and $p_\bullet$. Moreover, if the stated assertions are known for
all monomorphisms of functors and all $V$-linear categories $\cA$,
they follow for every epimorphism as well, by considering the
opposite transformation $f^o_\bullet:M'^o_\bullet\Rightarrow M^o_\bullet$
and invoking remark \ref{rem_V-lin-abelian-cat}(ii).

Thus, let $f_\bullet$ be a monomorphism, $\phi:j\to i$ any morphism
in $I$, and $a\in V$ that annihilates the morphism
$\Coker\,f_j\to\Coker\,f_i$ induced by $M'_\phi$. The latter means
that $\Img(aM'_\phi)\subset\Img\,f_i$; it follows that there exists
a morphism $g:M'_j\to M_i$ in $\cA$ such that $f_i\circ g=aM'_\phi$,
and therefore
$$
f_i\circ g\circ f_j=aM'_\phi\circ f_j=af_i\circ M_\phi=f_i\circ(aM_\phi)
$$
whence $g\circ f_j=aM_\phi$, since $f_i$ is a monomorphism. The
induced identities $Ff_i\circ Fg=aFM'_\phi$ and $Fg\circ Ff_j=aFM_\phi$
imply that $\Img(aFM'_\phi)\subset\Img\,Ff_i$ and
$aM_\phi(\Ker\,Ff_j)=0$. Hence, if $\Ker\,f_\bullet$ and
$\Coker\,f_\bullet$ are almost essentially zero (resp. dual
almost essentially zero) then the same holds for $\Ker\,Ff_\bullet$
and $\Coker\,Ff_\bullet$. Combining with (ii), we deduce that if
$\Ker\,f_\bullet$ and $\Coker\,f_\bullet$ are null (resp. dual null),
the same follows for $\Ker\,Ff_\bullet$ and $\Coker\,Ff_\bullet$.

(iii) is an immediate consequence of (i).

(iv): As usual, the assertion in case $F$ is left exact will
follow from the assertion for the case $F$ is right exact, by
considering $M^o_\bullet$ and $F^o$. Hence, suppose that $F$ is
right exact and $M_\bullet$ has the almost Mittag-Leffler
property. For a given subideal $\fm_0\subset\fm$ of finite
type and $i\in\Ob(I)$, choose a morphism $\phi:j\to i$ in
$I$ such that for every morphism $\psi:k\to j$ of $I$ we
have $\fm_0\cdot\Img\,M_\phi\subset\Img\,M_{\phi\circ\psi}$.
Let $C$ be the cokernel of the morphism
$\bar M_{\phi\circ\psi}:M_k\to\Img\,M_\phi$ of $\cA$ induced by
$M_{\phi\circ\psi}$; the assumption means that $\fm_0C=0$. Since
$F$ is right exact, $FC$ is the cokernel of
$F\bar M_{\phi\circ\psi}:FM_k\to F(\Img\,M_\phi)$, and obviously
$\fm_0FC=0$. By the same token, $M_\phi:M_j\to M_i$ factors
through epimorphisms $FM_j\to F(\Img\,M_\phi)\to\Img\,FM_\phi$;
we conclude that the cokernel of the morphism
$FM_k\to\Img\,FM_\phi$ induced by $FM_{\phi\circ\psi}$ is a
quotient of $FC$, and especially is annihilated by $\fm_0$.
The assertion follows.
\end{proof}

\begin{remark}\label{rem_dual-of-a.ML}
Keep the notation of definition \ref{def_almost-Mittag-Leff},
and let $M_\bullet:I\to\cA$ be any functor.

(i)\ \
By unwinding the definitions, we find that $M_\bullet$ has
the dual almost Mittag-Leffler property if and only if the
following holds. For every subideal $\fm_0\subset\fm$ of
finite type, and every $i\in\Ob(I)$ there exists a morphism
$\phi:i\to j$ in $I$ such that
$$
\fm_0\cdot\Ker\,M_{\psi\circ\phi}\subset\Ker\,M_\phi
\qquad
\text{for every morphism $\psi:j\to k$ in $I$}.
$$

(ii)\ \
The natural isomorphism $\bFun(I,\cA)\isom\bFun(I^o,\cA^o)^o$
of remark \ref{rem_opposite-Fun}(i) identifies
$(\bFun(I^o,\cA^o)_\mathrm{a.ess.0})^o$ with the full subcategory
of $\bFun(I,\cA)$ whose objects are the dual almost essentially
zero functors. Since the opposite of a Serre subcategory is a
Serre subcategory, lemma \ref{lem_a.ess.0-is-Serre}(i) implies
that the quotient
$$
\cC(I,\cA):=\bFun(I,\cA)/(\bFun(I^o,\cA^o)_\mathrm{a.ess.0})^o
$$
is again an abelian category, whose objects are the functors
$I\to\cA$.

(iii)\ \
The rest of lemma \ref{lem_a.ess.0-is-Serre} dualizes as well :
first, part (ii) of the lemma implies that if $\cA$ is
$\fm$-torsion-free and $M_\bullet$ is dual almost essentially
zero, then $\colim_IM_\bullet=0$. Next, if $I$ is filtered, and
$\cA$ is cocomplete and $\fm$-torsion-free, the functor
$\Colim_I$ of remark \ref{rem_Lims-and-Colims}(ii) factors
through the localization $\bFun(I,\cA)\to\cC(I,\cA)$ and a
functor
$$
\cC(I,\cA)\to\cA
\qquad
M_\bullet\mapsto\colim_IM_\bullet.
$$

(iv)\ \
Also remark \ref{rem_almost-Mittag-Leff}(ii,iii) admits
a valid dual : namely, suppose that $I$ is filtered and
$\cA$ is cocomplete. We set
$$
M^\triangledown_i:=\bigcup_{\phi:i\to j}\Ker\,M_\phi\subset M_i
\qquad
\text{for every $i\in\Ob(I)$}.
$$
Then, arguing as in \ref{rem_almost-Mittag-Leff}(ii) one
easily checks that $M_\phi$ restricts to a morphism
$$
M^\triangledown_\psi:M_i^\triangledown\to M^\triangledown_{i'}
\qquad
\text{for every morphism $\psi:i\to i'$ in $I$}
$$
whence an induced functor $M^\triangledown_\bullet:I\to\cA$,
and a direct inspection of the definitions shows that
\set\begin{equation}\label{eq_double-triangle}
(M^\triangledown_\bullet)^\triangledown_\bullet=M^\triangledown_\bullet.
\end{equation}
Notice moreover that every cocone $\tau:M_\bullet\Rightarrow c_L$
factors (uniquely) through the projection
$\pi_\bullet:M_\bullet\to M_\bullet/M^\triangledown_\bullet$. Combining
with \eqref{eq_double-triangle} we deduce that
$$
\colim_IM^\triangledown_\bullet=0
$$
and furthermore, $\pi_\bullet$ induces an isomorphism in $\cA$ :
$$
\colim_IM_\bullet\isom\colim_IM_\bullet/M^\triangledown_\bullet.
$$

(v)\ \
Under the assumptions of (iv), we deduce that $M_\bullet$
has the dual almost Mittag-Leffler property if and only if
for every subideal $\fm_0\subset\fm$ of finite type and every
$i\in\Ob(I)$ there exists a morphism $\phi:i\to j$ such that
$\fm_0\cdot M^\triangledown_i\subset\Ker\,M_\phi$. In the
case where $\cA=A\Mod$ for a given $(V,\fm)^a$-algebra
$A$, it follows that for every $i\in\Ob(I)$ the system
$(\Ker\,M_\phi~|~\phi:i\to j)$ is a Cauchy net in the
uniform space $\cI_A(M_i)$, and its unique limit is the
$A$-submodule $M^\triangledown_i$.

(vi)\ \
Suppose that $I$ is filtered, and is $\cA$ is cocomplete
and $\fm$-torsion-free. Then :

$\bullet$\ \
By dualizing proposition \ref{prop_a-Mittag-Leffler}(ii),
we see that the following conditions are equivalent :
\begin{enumerate}
\alphaenu
\item
$M_\bullet$ is dual almost essentially constant.
\item
$M_\bullet$ is isomorphic to a constant functor, in the category
$\cC(I,\cA)$ of remark \ref{rem_dual-of-a.ML}(ii).
\item
$M_\bullet$ has the dual almost Mittag-Leffler property, and for
every subideal $\fm_0\subset\fm$ of finite type, there exists
$i\in\Ob(I)$ such that $\fm_0$ annihilates the cokernel of the
morphism
\set\begin{equation}\label{eq_tutto-bloccato}
M_i/M^\triangledown_i\to M_j/M^\triangledown_j
\end{equation}
induced by $M_\phi$, for every morphism $\phi:i\to j$ in $I$
(notation of remark \ref{rem_dual-of-a.ML}(iv)).
\end{enumerate}

$\bullet$\ \
Moreover, the dual of proposition \ref{prop_a-Mittag-Leffler}(i)
states that if $M_\bullet$ has the dual almost Mittag-Leffler
property, then \eqref{eq_tutto-bloccato} is a monomorphism
for every morphism $\phi:i\to j$ in $I$, and the projections
$(M_i\to M_i/M^\triangledown_i~|~i\in\Ob(I))$ define an
isomorphism $M_\bullet\isom M_\bullet/M^\triangledown_\bullet$
in $\cC(I,\cA)$.
\end{remark}

\begin{lemma}\label{lem_lim-one-special-case}
In the situation of definition {\em\ref{def_almost-Mittag-Leff}},
let $A$ be any $(V,\fm)^a$-algebra, $M_\bullet:\N^o\to A\Mod$ a
functor whose transition morphisms $M_{m,n}:M_m\to M_n$ are
epimorphisms for every $m,n\in\N$, and set $L:=\lim_{n\in\N}M_n$.
Then we have :
\begin{enumerate}
\item
The universal cone $\tau_\bullet:c_L\Rightarrow M_\bullet$ is an
epimorphism of functors.
\item
$\lim^1_{n\in\N}M_n=0$.
\item
For every almost finitely generated (resp. almost finite
projective) $A$-module $N$, the cone $N\otimes_A\tau_\bullet$
induces an epimorphism (resp. an isomorphism)
$$
\phi_N:N\otimes_AL\to\lim_{n\in\N}\,N\otimes_AM_n.
$$
\item
More generally, if $N$ is an $A$-module and $a\in A_*$ such
that $a\cdot\one_N$ factors as a composition of $A$-linear
morphisms $N\xrightarrow{\alpha}N_0\xrightarrow{\beta}N$ with
$N_0$ almost finitely generated (resp. almost finite projective),
then $a\cdot\Coker\,\phi_N=0$ (resp. $a\cdot\Ker\,\phi_N=0$).
\end{enumerate}
\end{lemma}
\begin{proof}(i): Let $M_{\bullet!}:\N^o\to A_*\Mod$ be the
composition of $M_\bullet$ with the left adjoint
$(-)_!:A\Mod\to A_*\Mod$ of the localization $(-)^a$.
Set $T:=\lim_{n\in\N}M_{n!}$, and let
$\mu_\bullet:c_T\Rightarrow M_{\bullet!}$ be the universal cone;
then $L=T^a$ and the natural isomorphism
$M_\bullet\isom(M_{\bullet!})^a$ identifies $\tau_\bullet$ with
$\mu^a_\bullet$. Thus, it suffices to check that $\mu_\bullet$
is an epimorphism of functors. However, the transition
morphisms $(M_{m,n})_!:M_{m!}\to M_{n!}$ are surjective, so
the assertion is clear.

(ii): Recall that for every functor $M_\bullet:\N^o\to A\Mod$
(resp. $N_\bullet:\N^o\to A_*\Mod$) the right derived functor
$R\lim_{\N^o}M_\bullet$ (resp. $R\lim_{\N^o}N_\bullet$) is computed
by a natural two-term complex $\Pi(M_\bullet)$ (resp.
$\Pi(N_\bullet)$) : see \cite[\S3.5]{We}. With this notation,
we have $\Pi(M_\bullet)=\Pi(M_{\bullet!})^a$, since the
localization $(-)^a$ commutes with all products. However,
$\lim^1_{n\in\N}M_{n!}=0$, since the transition morphisms
$(M_{m,n})_!$ are surjective (\cite[Lemma 3.5.3]{We}). The
assertion follows.

(iii): Suppose first that $N$ is finitely generated; hence
there exists $r\in\N$ and a short exact sequence
$0\to K\to A^{\oplus r}\to N\to 0$ of $A$-modules. For every
$n\in\N$, let $C_n\subset A^{\oplus r}\otimes_AM_n$ be the image
of $K\otimes_AM_n$; we get an inverse system $(C_n~|~n\in\N)$
whose transition morphisms are also epimorphisms (details left
to the reader). By virtue of (ii), there follows an epimorphism
of $A$-modules
$$
\lim_{n\in\N}A^{\oplus r}\otimes_AM_n=A^{\oplus r}\otimes_AL
\xrightarrow{\psi}\lim_{n\in\N}N\otimes_AM_n.
$$
Lastly, it is easily seen that $\psi$ factors through $\phi_N$
and the induced morphism $A^{\oplus r}\otimes_AL\to N\otimes_AL$,
whence the contention in this case. Next, consider the functor
$$
T:A_*\Mod\to A_*\Mod
\qquad
X\mapsto(\Coker\,\phi_{X^a})_*
$$
that assigns to every $A$-linear morphism $f:N\to N'$ the
induced map
$(\Coker\,\phi_{N^a})_*\to(\Coker\,\phi_{N'^a})_*$ of $A_*$-modules.
It is easily seen that $T$ is $V$-homogeneous of degree $1$ (see
definition \ref{def_V-homogeneous-fct}), and the associated functor
$T^a:A\Mod\to A\Mod$ is given by the rule:
$N\mapsto\Coker\,\phi_N$ for every $A$-module $N$ (notation
of remark \ref{rem_V-homogeneous-funct}(ii)).
By the foregoing, $T^aN=0$ whenever $N$ is finitely generated;
however, on the one hand, for every suitable pair of infinite
cardinalities $\omega,\omega'$ the associated map
$\cM_{\omega,\omega'}(T^a):\cM_\omega(A)\to\cM_{\omega'}(A)$ is uniformly
continuous (proposition \ref{prop_ext-std-to-m-nonflat}(iii)), and
on the other hand, it is easily seen that the subset $\{0\}$ is
closed in the topology of the uniform space $\cM(A)$. Hence,
$T^aN=0$ for every almost finitely generated $A$-module, as stated.

Next, suppose that $N$ is almost finite projective; then
for every $a\in\fm$ there exists $r\in\N$ and $A$-linear
morphisms $N\xrightarrow{\alpha}A^{\oplus r}\xrightarrow{\beta}N$
whose composition is $a\cdot\one_N$ (\cite[Lemma 2.4.15]{Ga-Ra}).
There follows a commutative diagram :
$$
\xymatrix@C+20pt{
N\otimes_AL \ar[r]^-{\alpha\otimes_AL} \ar[d]_{\phi_N} &
A^{\oplus r}\otimes_AL \ar[r]^-{\beta\otimes_AL} \ar[d] &
N\otimes_AL \ar[d]^{\phi_N} \\
\lim_{n\in\N}N\otimes_AM_n \ar[r] &
\lim_{n\in\N}A^{\oplus r}\otimes_AM_n \ar[r] & \lim_{n\in\N}N\otimes_AM_n
}$$
whose central arrow is an isomorphism, and we already know that
$\phi_N$ is an epimorphism. By a simple diagram chase we deduce
that $a\cdot\Ker\,\phi_N=0$, and since $a$ is arbitrary, we
conclude that $\phi_N$ is an isomorphism in this case, as stated.

(iv): In light of (iii) we have :
$a\cdot\one_{TN}=T\beta\circ T\alpha=0$, {\em i.e.}
$a\cdot\Coker\,\phi_N=0$. If $N_0$ is almost finite projective,
consider the functor $T':A^a\Mod\to A^a\Mod$ such that
$T'X:=\Ker\,\phi_X$ for every $A^a$-module $X$ : by (iii) we
have $T'N_0=0$, whence $a\cdot\one_{T'N}=0$.
\end{proof}

\begin{proposition}\label{prop_last-ML-brick}
In the situation of definition {\em\ref{def_almost-Mittag-Leff}},
let $A$ be any $(V,\fm)^a$-algebra, $M_\bullet:I\to A\Mod$
any functor, and $\tau_\bullet:c_L\Rightarrow M_\bullet$ and
$\tau'_\bullet:M_\bullet\Rightarrow c_{L'}$ be respectively
a universal cone and a universal cocone. We have :

{\em(i)}\ \
Suppose that $I$ is cofiltered and consider the following
conditions :
\begin{enumerate}
\alphaenu
\item
The functor $\Coker\,\tau_\bullet$ is almost essentially zero.
\item
$M_\bullet$ is almost Mittag-Leffler.
\end{enumerate}
Then {\em(a)}$\Rightarrow${\em(b)}, and if there exists a
final functor $\N^o\to I$, then {\em(b)}$\Rightarrow${\em(a)}.

{\em(ii)}\ \
Suppose that is $I$ filtered. Then $M_\bullet$ has the dual almost
Mittag-Leffler property if and only if the functor
$\Ker\,\tau'_\bullet$ is dual almost essentially zero.

{\em(iii)}\ \
Suppose that $I=\N^o$, where $\N$ is endowed with its standard
total ordering. If $M_\bullet$ has the almost Mittag-Leffler
property, then $\lim_{n\in\N}^1M_n=0$.
\end{proposition}
\begin{proof}(iii): Suppose first that $M_\bullet$ is almost
essentially zero, and consider, for every finitely generated
subideal $\fm_0\subset\fm$, the short exact sequence of functors
$$
\Sigma_\bullet\qquad :\qquad
0\to M_\bullet[\fm_0]\to M_\bullet\to M_\bullet/M_\bullet[\fm_0]\to 0
$$
whose right-most term is an essentially zero inverse
system, by remark \ref{rem_almost-Mittag-Leff}(i)), and
therefore
$$
\lim_{n\in\N}M_n/M_n[\fm_0]=\lim_{n\in\N}{}^{\!1}M_n/M_n[\fm_0]=0.
$$
On the other hand, clearly $\fm_0$ annihilates
$\lim_{n\in\N}^1M_n[\fm_0]$. From the long exact cohomology
sequence associated with $\Sigma_\bullet$ it follows that
$\fm_0$ annihilates as well $\lim_{n\in\N}^1M_n$. Since
$\fm_0$ is arbitrary, the claim follows in this case.

Next, let $M_\bullet$ be an arbitrary functor with the almost
Mittag-Leffler property; for every $m,n\in\N$ with $m>n$, we
let $M_{m,n}:M_m\to M_n$ be the transition morphism. By
proposition \ref{prop_a-Mittag-Leffler}(i.a), the restriction
$M^\triangle_{p,n}:M_m^\triangle\to M^\triangle_n$ of $M_{m,n}$
is an epimorphism for every such $n$ and $p$. Let us then
consider the short exact sequence of inverse systems
$$
\Sigma'_\bullet \qquad :\qquad
0\to M^\triangle_\bullet\to M_\bullet\to M_\bullet/M^\triangle_\bullet\to 0.
$$
Proposition \ref{prop_a-Mittag-Leffler}(i.b) implies that
$M_\bullet/M^\triangle_\bullet$ is almost essentially
zero, so $\lim^1_{n\in\N}M_n/M^\triangle_n=0$, by the previous
case; also, $\lim^1_{n\in\N} M^\triangle_n=0$, by lemma
\ref{lem_lim-one-special-case}(ii). The assertion follows
then, by considering the long exact $R\lim$ sequence
associated with $\Sigma'_\bullet$.

(i.a)$\Rightarrow$(i.b): By assumption, for
every subideal $\fm_0\subset\fm$ of finite type and every
$i\in\Ob(I)$ there exists a morphism $\phi:j\to i$ in $I$
such that $\fm_0$ annihilates the image of the induced
morphism $\Coker\,\tau_j\to\Coker\,\tau_i$. The latter
means that $\fm_0\cdot\Img\,M_\phi\subset\Img\,\tau_i$;
but clearly $\Img\,\tau_i\subset\Img\,M_{\phi\circ\psi}$ for
every morphism $\psi:k\to j$ in $I$, whence the contention.

Next, suppose that there exists a final functor $\lambda:\N^o\to I$,
and that $M_\bullet$ has the almost Mittag-Leffler property;
notice that $\tau_i$ factors through a morphism
$\mu_i:L\to M^\triangle_i$ and the inclusion $M_i^\triangle\to M_i$,
for every $i\in\Ob(I)$. There follows a short exact sequence
of functors :
$$
\Coker\,\mu_\bullet\to\Coker\,\tau_\bullet\to
M_\bullet/M^\triangle_\bullet\to 0.
$$
By proposition \ref{prop_a-Mittag-Leffler}(i.b),
$M_\bullet/M^\triangle_\bullet$ is almost essentially
zero, hence $L$ also represents the limit of $M^\triangle_\bullet$,
and $\mu_\bullet$ is a universal cone; by lemma
\ref{lem_a.ess.0-is-Serre}(i) it then suffices to show that
$\Coker\,\mu_\bullet=0$.  Now,
$\mu_\bullet*\lambda:c_L\Rightarrow M^\triangle_\bullet\circ\lambda$
is still a universal cone (remark \ref{rem_fun-cofinal}(ii,iii)).
Since we also know that $M^\triangle_\phi$ is an epimorphism for
every morphism $\phi$ of $I$ (proposition
\ref{prop_a-Mittag-Leffler}(i.a)), we are then reduced to
checking that $\mu_{\lambda(n)}:L\to M^\triangle_{\lambda(n)}$
is an epimorphism for every $n\in\N$. The latter follows
from lemma \ref{lem_lim-one-special-case}(i).

(ii): Suppose that $\Ker\,\tau'_\bullet$ is dual almost
essentially zero; this means that for every subideal
$\fm_0\subset\fm$ of finite type, and every $i\in\Ob(I)$
there exists a morphism $\phi:i\to j$ in $I$ such that
$M_\phi(\fm_0\cdot\Ker\,\tau'_i)=0$. But we have
$\Ker\,M_{\psi\circ\phi}\subset\Ker\,\tau'_i$ for every
morphism $\psi:j\to k$ in $I$; whence
$\fm_0\cdot\Ker\,M_{\psi\circ\phi}\subset\Ker\,M_\phi$
for every such $\psi$. Hence $M_\bullet$ has the dual
Mittag-Leffler property, by remark \ref{rem_dual-of-a.ML}(i).
Conversely, suppose that $M_\bullet$ has the dual Mittag-Leffler
property; we know that the functor $M^\triangledown_\bullet$ is
almost essentially zero (remark \ref{rem_dual-of-a.ML}(vi))
hence it suffices to check that $M^\triangledown_i=\Ker\,\tau'_i$
for every $i\in\Ob(I)$. After replacing $I$ by $i/I$, we may
assume that $i$ is an initial object of $I$, whence an induced
cone $\mu_\bullet:c_{M_i}\Rightarrow M_\bullet$, and notice that
$M^\triangledown_i=\colim_I\Ker\,\mu_\bullet$. On the other hand,
$\colim_I\mu_\bullet=\tau'_i$; since filtered colimits in
$A\Mod$ are exact, the assertion follows.
\end{proof}

\sset\subsubsection{Pontryagin duality}\label{subsec_Pontryagin}
Recall that for every ring $R$ and every $R$-module $N$,
we have a natural $R$-module structure on
$$
N^\vee:=\Hom_\Z(N,\Q/\Z).
$$
Namely, for every $a\in R$ and every $\Z$-linear map
$\phi:N\to\Q/\Z$ we set $a\cdot\phi:=\phi\circ(a\one_N)$.
Especially, $R^\vee$ is naturally an $R$-module, and there
follows a natural $R$-linear isomorphism :
$$
\Hom_R(N,R^\vee)\isom N^\vee
$$
assigning to every $R$-linear map $\phi:N\to R^\vee$
the $\Z$-linear map $\phi':N\to\Q/\Z$ such that
$\phi'(n):=\phi(n)(1)$ for every $n\in N$.
Its inverse assigns to every $\Z$-linear map
$\psi:N\to\Q/\Z$ the homomorphism of $R$-modules
$\psi':N\to R^\vee$ that sends every $n\in N$ to
the $\Z$-linear map $\psi'(n):R\to\Q/\Z$ such that
$\psi'(n)(a):=\psi(an)$ for every $a\in R$.
Since $\Q/\Z$ is an injective $\Z$-module, there
results an exact functor
$$
R\Mod\to R\Mod^o
\qquad
N\mapsto N^\vee.
$$
Moreover, it is easily seen that the natural {\em biduality
$R$-linear map} is an injection :
$$
N\to N^{\vee\vee}
\qquad
\text{for every $R$-module $N$}.
$$
Let now $(V,\fm)$ be a basic setup, and $A$ a
$(V,\fm)^a$-algebra; we consider the associated $A$-module
$$
J_A:=((A_*)^\vee)^a
$$
and from the foregoing it is easily seen that the
{\em Pontryagin duality} functor
$$
(-)^\vee:A\Mod\to A\Mod^o
\qquad
M\mapsto M^\vee:=\Alhom_A(M,J_A)
$$
is again exact and $V$-linear, where $\Alhom_A(-,-)$ denotes
the bifunctor of {\em almost homomomorphisms} defined as in
\cite[\S2.2.11]{Ga-Ra}; indeed, we have a natural isomorphism
of $A$-modules
$$
M^\vee\isom((M_*)^\vee)^a
\qquad
\text{for every $A$-module $M$}
$$
(details left to the reader). Moreover, it follows that the
biduality morphism $M\to M^{\vee\vee}$ is a monomorphism of
$A$-modules; especially, $(-)^\vee$ is a conservative functor.

\begin{remark}\label{rem_Pontryagin-for-functors}
(i)\ \
For every small category $I$, we get a natural extension
of Pontryagin duality to functors $I\to A\Mod$. Namely,
we have a conservative exact functor :
$$
\bFun(I,A\Mod)^o\to\bFun(I^o,A\Mod)
\qquad
M_\bullet\mapsto M^\vee_\bullet:=(-)^\vee\circ M_\bullet.
$$

(ii)\ \
For every $A$-module $M$, and every $A$-submodule $N\subset M$,
let $i_N:N\to M$ the inclusion; it is easily seen that the
induced map
$$
\cI_A(M)\to\cI_A(M^\vee)
\qquad
N\mapsto\Ker(i^\vee_N:M^\vee\to N^\vee)
$$
is uniformly continuous, for the uniform structures
defined as in \cite[Def.2.3.1(i)]{Ga-Ra}.

(iii)\ \
Then, if $I$ is cofiltered, and $M_\bullet:I\to A\Mod$
is any functor with the almost Mittag-Leffler property,
combining with remarks \ref{rem_almost-Mittag-Leff}(iii)
and \ref{rem_dual-of-a.ML}(v) we deduce that the $A$-linear
epimorphism $i^\vee_{M^\triangle_i}:M_i^\vee\to(M^\triangle_i)^\vee$
factors through an isomorphism of $A$-modules
$$
M_i^\vee/(M^\vee_\bullet)^\triangledown_i\isom(M^\triangle_i)^\vee
\qquad
\text{for every $i\in\Ob(I)$}.
$$
\end{remark}

\begin{lemma}\label{lem_a.ess.cst-and-duals}
With the notation of remark
{\em\ref{rem_Pontryagin-for-functors}}, for every
functor $M_\bullet:I\to A\Mod$ the following holds :

{\em (i)}\ \
$M_\bullet$ is essentially zero (resp. almost essentially zero)
if and only if $M^\vee_\bullet$ is dual essentially zero (resp.
dual almost essentially zero).

{\em (ii)}\ \
Suppose that $I$ is cofiltered. Then $M_\bullet$ is almost
essentially constant (resp. has the almost Mittag-Leffler
property) if and only if $M^\vee_\bullet$ is dual almost essentially
constant (resp. has the dual almost Mittag-Leffler property).
\end{lemma}
\begin{proof}(i): By direct inspection we see that if
$M_\bullet$ is essentially zero (resp. almost essentially
zero), then $M^\vee$ is dual essentially zero (resp.
dual almost essentially zero). Likewise, if $M_\bullet$
is dual essentially zero (resp. dual almost essentially
zero), then $M^\vee$ is essentially zero (resp. almost
essentially zero). From this, the converse assertions
follow easily, since the biduality morphism
$M_\bullet\to M^{\vee\vee}_\bullet$ is a monomorphism.

(ii): For a given pair of morphisms
$k\xrightarrow{\psi}j\xrightarrow{\phi}i$ of $I$, set
$N:=\Img M_\phi/\Img\,M_{\phi\circ\psi}$. Since Pontryagin
duality is exact and conservative, for every subideal
$\fm_0\subset\fm$ we have :
$$
\fm_0N=0\Leftrightarrow\fm_0N^\vee=0\Leftrightarrow
\fm_0\cdot\Ker(\Img\,M^\vee_{\phi\circ\psi}\to M^\vee_\phi)=0
\Leftrightarrow
\fm_0\Ker\,M^\vee_{\phi\circ\psi}\subset\Ker M^\vee_\phi.
$$
Taking into account remark \ref{rem_dual-of-a.ML}(i),
we deduce that $M_\bullet$ has the almost Mittag-Leffler
property if and only if $M^\vee_\bullet$ has the dual almost
Mittag-Leffler property.

Next, proposition \ref{prop_apply-F-V-linear}(ii) implies
that if $M_\bullet$ is almost essentially constant, then
$M^\vee_\bullet$ is dual almost essentially constant.
Conversely, if $M^\vee_\bullet$ is dual almost essentially
constant, then $M^\vee_\bullet$ has the dual almost
Mittag-Leffler property, and for every subideal $\fm_0\subset\fm$
of finite type there exists $i\in\Ob(I)$ with
\set\begin{equation}\label{eq_stupid-fly}
\fm_0\cdot\Coker(M^\vee_i/(M^\vee_\bullet)^\triangledown_i\to
M^\vee_j/(M^\vee_\bullet)^\triangledown_j)=0
\qquad
\text{for every morphism $i\xrightarrow{\phi}j$}
\end{equation}
(remark \ref{rem_dual-of-a.ML}(vi)). By the foregoing, we
deduce already that $M_\bullet$ has the almost Mittag-Leffler
property. Moreover, by remark
\ref{rem_Pontryagin-for-functors}(iii), condition
\eqref{eq_stupid-fly} is equivalent to :
$$
\fm_0\cdot\Coker((M^\triangle_i)^\vee\to(M^\triangle_j)^\vee)=0
\qquad
\text{for every morphism $i\xrightarrow{\phi}j$}.
$$
But $\Coker((M^\triangle_i)^\vee\to(M^\triangle_j)^\vee)=
(\Ker\,M^\triangle_\phi)^\vee$, and in light of proposition
\ref{prop_a-Mittag-Leffler}(ii), we deduce that $M_\bullet$
is almost essentially constant.
\end{proof}

\begin{proposition}\label{prop_clear-from-def}
Let $I$ be a small cofiltered category, and
$$
0\to M'_\bullet\to M_\bullet\to M''_\bullet\to 0
$$
a short exact sequence of functors $I\to A\Mod$. We have:

{\em (i)}\ \
If $M'_\bullet$ and $M''_\bullet$ are almost essentially constant
(resp. have the almost Mittag-Leffler property), then the same
holds for $M_\bullet$.

{\em (ii)}\ \
If $M_\bullet$ has the almost Mittag-Leffler property, the same
holds for $M''_\bullet$.

{\em (iii)}\ \
If $M_\bullet$ has the almost Mittag-Leffler property and
$M''_\bullet$ is almost essentially constant, then $M'_\bullet$
has the almost Mittag-Leffler property.

{\em (iv)}\ \
If $f_\bullet:X_\bullet\to Y_\bullet$ is a natural transformation
between almost essentially constant functors, then $\Ker\,f_\bullet$
and $\Coker\,f_\bullet$ are almost essentially constant.
\end{proposition}
\begin{proof}(ii) is clear from the definitions.

(iv): By proposition \ref{prop_a-Mittag-Leffler}(ii),
the universal cones $\tau^X_\bullet:c_X\Rightarrow X_\bullet$
and $\tau^Y_\bullet:c_X\Rightarrow Y_\bullet$ are isomorphisms
in $\cL(I,A\Mod)$. Thus, there exists a unique morphism
$f:X\to Y$ of $A$-modules that makes commute the diagram
$$
\xymatrix{ c_X \ar[r]^-{c_f} \ar[d]_{\tau^X_\bullet} &
c_Y \ar[d]^{\tau^Y_\bullet} \\
X_\bullet \ar[r]^-{f_\bullet} & Y_\bullet
}$$
whose vertical arrows are isomorphism in $\cL(I,A\Mod)$.
Set $K:=\Ker\,f$ and $C:=\Coker\,f$; there follow induced
isomorphisms $c_K\isom\Ker\,f_\bullet$ and
$c_C\isom\Coker\,f_\bullet$ in $\cL(I,A\Mod)$, whence the
contention, in light of proposition
\ref{prop_a-Mittag-Leffler}(ii).

(i): Let $L,L'$ and $L''$ be the colimits of the functors
$M^\vee_\bullet,M'^\vee_\bullet$ and respectively $M''^\vee_\bullet$;
since filtered colimits are exact in the category $A\Mod$,
we get a commutative diagram of functors with exact rows :
$$
\xymatrix{ 0 \ar[r] & M''^\vee_\bullet \ar[r] \ar[d] &
M^\vee_\bullet \ar[r] \ar[d] & M'^\vee_\bullet \ar[r] \ar[d] & 0 \\
0 \ar[r] & c_{L''} \ar[r] & c_L \ar[r] & c_{L'} \ar[r] & 0
}$$
whose vertical arrows are the universal cocones. Consider
first the case where $M'_\bullet$ and $M''_\bullet$ are almost
essentially constant; then $M'^\vee_\bullet$ and $M''^\vee_\bullet$
are dual almost essentially constant (lemma
\ref{lem_a.ess.cst-and-duals}(ii)), and arguing as in the
proof of (iv) we deduce that the first and third vertical
arrows are isomorphisms in the quotient category $\cC(I,A\Mod)$,
hence the same holds for the middle one, by the 5-lemma.
Then $M^\bullet_\bullet$ is dual almost essentially constant,
by remark \ref{rem_dual-of-a.ML}(vi), hence $M_\bullet$ is
almost essentially constant, again by lemma
\ref{lem_a.ess.cst-and-duals}(ii).

Next, if $M'_\bullet$ and $M''_\bullet$ have the almost
Mittag-Leffler property, then $M'^\vee_\bullet$ and $M''^\vee_\bullet$
have the dual almost Mittag-Leffler property (lemma
\ref{lem_a.ess.cst-and-duals}(ii)), hence the kernels
of the first and third vertical arrows are dual almost
essentially zero (proposition \ref{prop_last-ML-brick}(ii)).
Then the same holds for the kernel of the middle vertical
arrow, so $M^\vee_\bullet$ has the dual almost Mittag-Leffler
property, again by proposition \ref{prop_last-ML-brick}(ii);
finally, $M_\bullet$ has the almost Mittag-Leffler property,
by lemma \ref{lem_a.ess.cst-and-duals}(ii).

(iii): We consider again the induced commutative ladder
with exact rows. Arguing as in the proof of (ii), we see
that the kernel of the middle vertical arrow and the cokernel
of the left-most vertical arrow are both dual almost
essentially zero, hence the same holds for the kernel of
the right-most vertical arrow, by the snake lemma. Then
$M'^\vee_\bullet$ has the dual almost Mittag-Leffler property,
by proposition \ref{prop_last-ML-brick}(ii), and finally,
$M'_\bullet$ has the almost Mittag-Leffler property, by
lemma \ref{lem_a.ess.cst-and-duals}(ii).
\end{proof}

\sset\subsubsection{}\label{prop_inverse-syst-of-AlgMods}
Let us now consider a basic setup $(V,\fm)$, a (small)
cofiltered category $I$, and a functor
$$
I\to(V,\fm)^a\tdu\AlgMod
\qquad
i\mapsto(A_i,M_i)
$$
where $(V,\fm)^a\tdu\AlgMod$ is defined as in
\cite[Def.2.5.22(ii)]{Ga-Ra}. Hence, for every $i\in\Ob(I)$,
the pair $(A_i,M_i)$ consists of a $(V,\fm)^a$-algebra $A_i$,
and an $A_i$-module $M_i$. To each morphism $\phi:j\to i$ of
$I$, the functor assigns a pair
$(A_\phi,g_\phi):(A_j,M_j)\to(A_i,M_i)$, where $A_\phi:A_j\to A_i$
is a morphism of $(V,\fm)^a$-algebras, and
$g_\phi:A_i\otimes_{A_j}M_j\to M_i$ is an $A_i$-linear morphism.
We get two obvious induced functors
$$
A_\bullet:I\to(V,\fm)^a\Alg
\qquad
M_\bullet:I\to(V,\fm)^a\Mod.
$$
Namely, for every morphism $\phi$ as in the foregoing we
let $M_\phi:M_j\to M_i$ be the composition of $g_\phi$ with
the natural $A_j$-linear morphism $M_j\to A_i\otimes_{A_j}M_j$.
Set also
$$
A:=\lim_IA_\bullet
\qquad\text{and}\qquad
M:=\lim_IM_\bullet.
$$
Notice that $A_\bullet$ (resp. $M_\bullet$) can also be viewed as
a functor with values in $A\Alg$ (resp. $A\Mod$).

\begin{proposition}\label{prop_return-of-AlMod}
With the notation of \eqref{prop_inverse-syst-of-AlgMods},
the following holds :

{\em (i)}\ \
If $A_\bullet:I\to A\Mod$ has the almost Mittag-Leffler property,
and $g_\phi$ is an epimorphism for every morphism $\phi$ of $I$,
then $M_\bullet:I\to A\Mod$ has the almost Mittag-Leffler
property.

{\em (ii)}\ \
Suppose that $g_\phi$ is an isomorphism for every morphism
$\phi$ of $I$. Then, if\/ $A_\bullet:I\to A\Mod$ is almost
essentially constant (resp. is a Cauchy functor), the same
holds for $M_\bullet$.

{\em (iii)}\ \
Moreover, under the assumptions of {\em (i)} (resp. {\em (ii)}),
the natural morphism
$$
\alpha_i:A_i\otimes_AM_i^\triangle\to M_i
\qquad
\text{(resp. $\beta_i:A_i\otimes_AM\to M_i$\ )}
$$
is an epimorphism (resp. an isomorphism) for every $i\in\Ob(I)$.

{\em (iv)}\ \
Lastly, under the assumptions of {\em (i)}, if there exists a
final functor $\N^o\to I$, then $\beta_i$ is an epimorphism for
every $i\in\Ob(I)$.
\end{proposition}
\begin{proof} We fix universal cones
$\tau^A_\bullet:c_A\Rightarrow A_\bullet$ and
$\tau^M_\bullet:c_M\Rightarrow M_\bullet$.

\begin{claim}\label{cl_claimification}
Assertions (ii) and (iii) hold if $A_\bullet$ is a Cauchy functor and
$g_\phi$ is an isomorphism for every morphism $\phi$ of $I$.
\end{claim}
\begin{pfclaim} Indeed, in this case for every subideal
$\fm_0\subset\fm$ of finite type we find $i\in\Ob(I)$ such
that for every morphism $\phi:j\to i$ we have
$\fm_0\cdot\Ker\,A_\phi=\fm_0\cdot\Coker\,A_\phi=0$. By virtue
of remark \ref{rem_V-lin-abelian-cat}(vi), it follows that
$\fm_0^2$ annihilates the kernel and cokernel of
$A_\phi\otimes_{A_j}M_j$, for every such $\phi$. Since
$g_\phi$ is an isomorphism, it follows that
$\fm_0^2\cdot\Ker\,M_\phi=\fm_0^2\cdot\Coker\,M_\phi=0$, for
every such $\phi$. This shows that $M_\bullet$ is a Cauchy
functor.

Next, notice that $\Ker\,\tau^A_\bullet$ and
$\Coker\,\tau^A_\bullet$ are also Cauchy functors (remark
\ref{rem_Cauchy-functors}(ii)), and moreover they are
almost essentially zero, since $A_\bullet$ is almost essentially
constant (proposition \ref{prop_a-Mittag-Leffler}(ii,iii)).
Hence, $\Ker\,\tau^A_\bullet$ and $\Coker\,\tau^A_\bullet$ are
null. The same argument shows that $\Ker\,\tau^M_\bullet$
and $\Coker\,\tau^M_\bullet$ are null, since $M_\bullet$
is a Cauchy functor. Moreover, the kernel and
cokernel of $\tau^A_\bullet\otimes_AM$ are null
as well, since the functor $-\otimes_AM:A\Mod\to A\Mod$ is
$V$-linear (proposition \ref{prop_apply-F-V-linear}(i)).
Since we have :
$$
\beta_i\circ(\tau^A_i\otimes_AM)=\tau^M_i
\qquad
\text{for every $i\in\Ob(I)$}
$$
and as the null functors form a full Serre subcategory
of $\bFun(I,A\Mod)$ (remark \ref{rem_almost-Mittag-Leff}(v)), it
follows that also $\Ker\,\beta_\bullet$ and $\Coker\,\beta_\bullet$
are null functors. Now, let $\fm_0\subset\fm$ be a subideal
of finite type, and $i\in\Ob(I)$; since $I$ is cofiltered, we deduce
that there exists a morphism $\phi:j\to i$ such that
$\fm_0\cdot\Ker\,\beta_j=\fm_0\cdot\Coker\,\beta_j=0$.
Since the functor $A_i\otimes_{A_j}-:A_j\Mod\to A_i\Mod$ is
$V$-linear, it follows that $\fm_0^2$ annihilates the kernel
and cokernel of
$A_i\otimes_{A_j}\beta_j:A_i\otimes_AM\to A_i\otimes_{A_j}M_j$
(remark \ref{rem_V-lin-abelian-cat}(vi)). But notice that
$$
g_\phi\circ(A_i\otimes_{A_j}\beta_j)=\beta_i.
$$
So finally $\fm_0^2\cdot\Ker\,\beta_i=\fm_0^2\cdot\Coker\,\beta_i=0$.
Since $\fm_0$ is arbitrary, we conclude that $\beta_i$ is an
isomorphism, as stated.
\end{pfclaim}

We can now complete the proof of (ii) and (iii), for the case
where $A_\bullet$ is almost essentially constant, and $g_\phi$
is still an isomorphism for every morphism $\phi$ of $I$.
To this aim, set
$$
\bar A_\phi:=\Img\,A_\phi
\quad\text{and}\quad
\bar M_\phi:=\bar A_\phi\otimes_{A_j}M_j
\qquad
\text{for every morphism $\phi:j\to i$ of $I$}.
$$
Let $j\xrightarrow{\phi}i$ and $j'\xrightarrow{\phi'}i'$ be
two object of $\sMorph(I)$, and $(\psi,\psi'):\phi'\to\phi$
a morphism of $\sMorph(I)$; by definition, $\psi:j'\to j$ and
$\psi':i'\to i$ are morphisms in $I$ such that
$\psi'\circ\phi'=\phi\circ\psi$, and then it is clear that
$A_{\psi'}$ restricts to a morphism of $(V,\fm)^a$-algebras
$$
\bar A_{(\psi,\psi')}:\bar A_{\phi'}\to\bar A_\phi.
$$
Then, $M_\psi:M_{j'}\to M_j$ induces a $\bar A_\phi$-linear morphism
$$
\bar g_{(\psi,\psi')}:\bar A_\phi\otimes_{\bar A_{\phi'}}\bar M_{\phi'}
\isom\bar A_\phi\otimes_{A_{j'}}M_{j'}\to\bar M_\phi
$$
and it is easily seen that the rules :
$$
\phi\mapsto(\bar A_\phi,\bar M_\phi)
\qquad\text{and}\qquad
(\phi'\xrightarrow{(\psi,\psi')}\phi)\mapsto
(\bar A_{(\psi,\psi')},\bar g_{(\psi,\psi')})
$$
for every $\phi\in\Ob(\sMorph(I))$ and every morphism
$(\psi,\psi')$ of $\sMorph(I)$, define a functor
$$
\sMorph(I)\to(V,\fm)^a\tdu\AlgMod
$$
whence, as in \eqref{prop_inverse-syst-of-AlgMods}, two
functors
$$
\bar A_\bullet:\sMorph(I)\to A\Alg
\qquad
\bar M_\bullet:\sMorph(I)\to A\Mod
$$
and notice also that :
\set\begin{equation}\label{eq_dialetto}
\bar A{}^\triangle_\phi=A^\triangle_i
\qquad
\text{for every morphism $\phi:j\to i$ in $I$}.
\end{equation}
Next, let $\cI$ be the category whose objects are the
pairs $(j\xrightarrow{\phi}i,\fm_0)$ where $\phi$ is
a morphism of $I$ and $\fm_0\subset\fm$ is a subideal
of finite type such that :
$$
\fm_0\bar A_\phi\subset\bar A{}^\triangle_\phi
$$
(notation of remark \ref{rem_almost-Mittag-Leff}(ii)).
For every $(\phi,\fm_0),(\phi',\fm'_0)\in\Ob(\cI)$, the
set of morphisms $(\phi',\fm'_0)\to(\phi,\fm_0)$ in $\cI$
is empty if $\fm_0\not\subset\fm'_0$, and otherwise
consists of the morphisms $\phi'\to\phi$ in $\sMorph(I)$.
The composition law is then the same as that of $\sMorph(I)$,
so that the projection $(\phi,\fm_0)\mapsto\phi$ and the
rule : $i\mapsto\one_i$ for every $i\in\Ob(I)$ define functors
$$
\cI\xrightarrow{\pi}\sMorph(I)\xleftarrow{\delta}I.
$$

\begin{claim}\label{cl_pi-and-delta-final}
The functors $\pi$ and $\delta$ are final.
\end{claim}
\begin{pfclaim} It is easily seen that both $\cI$ and
$\sMorph(I)$ are cofiltered, so we can apply the (dual of the)
criterion of lemma \ref{lem_filtered-final}(i). Condition (a)
of the lemma is trivially verified by $\pi$. To check condition
(b) for $\pi$, consider any object $\phi:j\to i$ of $\sMorph(I)$,
any object $(\phi':j'\to i',\fm_0)$ of $\cI$, and a pair of
morphisms $(\psi,\psi'),(\mu,\mu'):\phi'\to\phi$. Since $I$
is cofiltered, we find morphisms $\nu:j''\to j'$ and
$\nu':i''\to i'$ in $I$ such that $\psi\circ\nu=\mu\circ\nu$ and
$\psi'\circ\nu'=\mu'\circ\nu'$. Then we also find $k\in\Ob(I)$
with morphisms $\lambda:k\to j''$ and $\lambda':k\to i''$
such that $\nu'\circ\lambda'=\phi'\circ\nu\circ\lambda$.
Set $\phi'':=\nu'\circ\lambda'$; we deduce a morphism
$(\nu\circ\lambda,\one_{i'}):\phi''\to\phi'$ in $\sMorph(I)$
with $(\psi,\psi')\circ(\nu\circ\lambda,\one_{i'})=
(\mu,\mu')\circ(\nu\circ\lambda,\one_{i'})$. From
\eqref{eq_dialetto} we get
$\bar A{}^\triangle_{\phi''}=\bar A{}^\triangle_{\phi'}$.
By assumption we have
$\fm_0\cdot\bar A_{\phi'}\subset\bar A{}^\triangle_{\phi'}$,
and clearly $\bar A_{\phi''}\subset\bar A_{\phi'}$. Hence
$\fm_0\cdot\bar A_{\phi''}\subset\bar A{}^\triangle_{\phi''}$, so
$(\phi'',\fm_0)\in\Ob(\cI)$, and $(\nu\circ\lambda,\one_{i'})$
is a morphism $(\phi'',\fm_0)\to(\phi',\fm_0)$ in $\cI$
such that $(\psi,\psi')\circ\pi(\nu\circ\lambda,\one_{i'})=
(\mu,\mu')\circ\pi(\nu\circ\lambda,\one_{i'})$, as required.
The verification of conditions (a) and (b) for the functor
$\delta$ is easy, and shall be left to the reader.
\end{pfclaim}

We consider then the induced functor :
$$
\bar A_\bullet\circ\pi:\cI\to A\Mod
\qquad
(j\xrightarrow{\phi}i,\fm_0)\mapsto\bar A_\phi.
$$
Notice that $\Img\,\tau^A_i\subset\bar A_\phi$ for every
morphism $\phi:j\to i$ of $I$; we deduce a cone
$$
\bar\tau{}^A_\bullet:c_A\Rightarrow\bar A_\bullet
\qquad
(j\xrightarrow{\phi}i)\mapsto(\tau^A_i:A\to\bar A_\phi)
$$
such that $\bar\tau{}^A_\bullet*\delta=\tau^A_\bullet$.
In light of claim \ref{cl_pi-and-delta-final}, it follows
that $\bar\tau{}^A_\bullet$ is a universal cone, and by the
same token, the same holds for the cone
$\bar\tau{}^A_\bullet*\pi:c_A\Rightarrow\bar A_\bullet\circ\pi$.

\begin{claim}\label{cl_Cauchification}
$\bar A_\bullet\circ\pi$ is a Cauchy functor.
\end{claim}
\begin{pfclaim} By proposition \ref{prop_a-Mittag-Leffler},
the functor $A^\triangle_\bullet:I\to A\Mod$ is Cauchy and
$A_\bullet$ has the almost Mittag-Leffler property, so
$A^\triangle_\phi:A^\triangle_j\to A^\triangle_i$ is an epimorphism
for every morphism $\phi:j\to i$ in $I$. Hence, for every
subideal $\fm_0\subset\fm$ of finite type there exists
$i\in\Ob(I)$ with $\fm_0\cdot\Ker\,A^\triangle_\phi=0$ for every
such $\phi$, and moreover there exists $\phi:j\to i$ with
$\fm_0\bar A_\phi\subset A^\triangle_i$. For such $\phi$,
we have $(\phi,\fm_0)\in\Ob(\cI)$; now, let
$(\psi,\psi'):(\phi':j'\to i',\fm'_0)\to(\phi,\fm_0)$
be any morphism in $\cI$; this means that $\fm_0\subset\fm'_0$
and $\fm'_0\bar A_{\phi'}\subset\bar A{}^\triangle_{\phi'}=A^\triangle_{i'}$.
Also, $\bar A_{(\psi,\psi')}:\bar A_{\phi'}\to\bar A_\phi$ is
the restriction of $A_{\psi'}:A_{i'}\to A_i$. Therefore
$$
\fm_0\cdot\Ker\,\bar A_{(\psi,\psi')}\subset
A^\triangle_{i'}\cap\Ker\,A_{\psi'}\subset\Ker\,A^\triangle_{\psi'}
$$
whence $\fm_0^2\cdot\Ker\,\bar A_{(\psi,\psi')}=0$. Lastly,
$\Coker\,\bar A_{(\psi,\psi')}$ is a quotient of
$\bar A_\phi/A_{\psi'}(A^\triangle_{i'})=\bar A_\phi/A^\triangle_i$,
which is annihilated by $\fm_0$, so
$\fm_0\cdot\Coker\,\bar A_{(\psi,\psi')}=0$.
\end{pfclaim}

Since the morphisms $g_\phi$ are isomorphisms for every
morphism $\phi$ of $I$, a simple inspection shows that
$\bar g_{(\psi,\psi')}$ is an isomorphism for every morphism
$(\psi,\psi')$ of $\sMorph(I)$. Together with claims
\ref{cl_Cauchification}, and \ref{cl_claimification},
we deduce that assertion (iii) holds for the functors
$\bar A_\bullet\circ\pi$ and $\bar M_\bullet\circ\pi$, and
moreover $\bar M_\bullet\circ\pi$ is a Cauchy functor.

However, let $\lambda:I\to\cI$ be the functor given by the
rules : $i\mapsto(\one_i,0)$ and $\phi\mapsto(\phi,\phi)$
for every $i\in\Ob(I)$ and every morphism $\phi$ of $I$. It
is easily seen that $\bar A_\bullet\circ\pi\circ\lambda=A_\bullet$
and $\bar M_\bullet\circ\pi\circ\lambda=M_\bullet$. Then assertion
(iii) for $A_\bullet$ and $M_\bullet$ follows straightforwardly.
Especially, the functor $M_\bullet$ is isomorphic to the
functor $A_\bullet\otimes_AM$. Since the functor
$-\otimes_AM:A\Mod\to A\Mod$ is $V$-linear, proposition
\ref{prop_apply-F-V-linear}(iii) implies then that $M_\bullet$
is almost essentially constant, and the proof of (iii) is
completed, under the assumptions of (ii).

Lastly, suppose that the assumptions of (i) hold; hence
$g_\phi$ is an epimorphism for every morphism $\phi:j\to i$
in $I$, and it follows easily that the induced morphism
$A_i\otimes_{A_j}\Img\,M_\phi\to M_i$ is an epimorphism as
well, for every such $\phi$. On the other hand, by assumption
for every $i\in\Ob(I)$ and every subideal $\fm_0\subset\fm$
of finite type there exists $\phi:j\to i$ such that
$\fm_0\Img\,M_\phi\subset M^\triangle_i$. Hence,
$\fm_0M_i\subset\Img\,\alpha_i$. Since $\fm_0$ is arbitrary,
we conclude that $\alpha_i$ is an epimorphism.

(i): By assumption, for every $i\in\Ob(I)$ and every
subideal $\fm_0\subset\fm$ there exists a morphism
$\phi:j\to i$ in $I$ such that for every morphism
$\psi:k\to j$ in $I$ we have
$\fm_0\cdot\Img\,A_{\phi\circ\psi}\subset\Img\,A_\phi$.
It then follows that for every such $\fm_0$, $\phi$
and $\psi$ we have
$$
\fm_0\cdot\Img(A_\phi\otimes_{A_k}M_k)\subset
\Img(A_{\phi\circ\psi}\otimes_{A_k}M_k).
$$
Since $g_\psi$ is an epimorphism and
$g_{\phi\circ\psi}\circ(A_{\phi\circ\psi}\otimes_{A_k}M_k)=
M_{\phi\circ\psi}$, we deduce :
$$
\fm_0\cdot\Img\,M_\phi=\fm_0\cdot\Img(M_\phi\circ g_\psi)=
\fm_0\cdot\Img(g_{\phi\circ\psi}\circ(A_\phi\otimes_{A_k}M_k))
\subset\Img\,M_{\phi\circ\psi}
$$
which shows that $M_\bullet$ has the almost Mittag-Leffler
property, as stated.

(iv): Notice that $\tau^M_i:M\to M_i$ factors through a
morphism $\mu_i:M\to M^\triangle_i$ and the inclusion
$M^\triangle_i\to M_i$, for every $i\in\Ob(I)$. Since
$M_\bullet/M^\triangle_\bullet$ is almost essentially zero
(proposition \ref{prop_a-Mittag-Leffler}(i.b)), it
follows that $\mu_\bullet:c_M\Rightarrow M^\triangle_\bullet$
is a universal cone. If $\lambda:\N^o\to I$ is a final
functor, the cone $\mu_\bullet*\lambda$ is still universal
(remark \ref{rem_fun-cofinal}(ii,iii)), and then
$\mu_{\lambda(n)}:M\to M^\triangle_{\lambda(n)}$ is an epimorphism
for every $n\in\N$ (lemma \ref{lem_lim-one-special-case}(i)).
Since $M^\triangle_\phi$ is an epimorphism for every $i\in\Ob(I)$
(proposition \ref{prop_a-Mittag-Leffler}(i.a)), we conclude
that $\mu_i:M\to M^\triangle_i$ is an epimorphism for every
$i\in\Ob(I)$. Since we already know that $\alpha_i$ is an
epimorphism, it follows that the same holds for $\beta_i$,
for every such $i$.
\end{proof}

\sset\subsubsection{}\label{subsec_Aretha-RIP}
Keep the notation of \eqref{prop_inverse-syst-of-AlgMods},
and consider a second functor
$$
I\to(V,\fm)^a\tdu\AlgMod
\qquad
i\mapsto(A'_i,N_i).
$$
as well as the corresponding induced functors
$$
A'_\bullet:I\to(V,\fm)^a\Alg
\qquad
N_\bullet:I\to(V,\fm)^a\Mod
$$
and suppose that $A_\bullet=A'_\bullet$. Let also $N$ be the
limit of $N_\bullet$. Then we get an induced functor
$$
I\to(V,\fm)^a\tdu\AlgMod
\qquad
i\mapsto(A_i,M_i\otimes_{A_i}N_i)
$$
and the corresponding induced functor
$$
M_\bullet\otimes_{A_\bullet}N_\bullet:I\to(V,\fm)^a\Mod
\qquad
i\mapsto M_i\otimes_{A_i}N_i.
$$

\begin{proposition}\label{prop_Aretha-RIP}
In the situation of \eqref{subsec_Aretha-RIP}, suppose that
$A_\bullet$, $M_\bullet$ and $N_\bullet$ are almost essentially
constant functors with values in $A\Mod$. Then the same holds
for $M_\bullet\otimes_{A_\bullet}N_\bullet$, and $M\otimes_AN$
represents the limit of $M_\bullet\otimes_{A_\bullet}N_\bullet$.
\end{proposition}
\begin{proof} For every $i\in\Ob(I)$, the $V^a$-linear
map
$$
\phi_i:M_i\otimes_{V^a}A_i\otimes_{V^a}N_i\to M_i\otimes_{V^a}N_i
\qquad
x\otimes a\otimes y\mapsto ax\otimes y-x\otimes ay
$$
yields a natural identification :
$\Coker\,\phi_i\isom M_i\otimes_{A_i}N_i$. In view of proposition
\ref{prop_clear-from-def}(iv), we are then easily reduced to
showing that the two functors :
$$
I\to(V,\fm)^a\Mod
\qquad
i\mapsto M_i\otimes_{V^a}A_i\otimes_{V^a}N_i
\qquad
i\mapsto M_i\otimes_{V^a}N_i
$$
are both almost essentially constant, and their limits
are represented by $M\otimes_{V^a}A\otimes_{V^a}N$ and
respectively $M\otimes_{V^a}N$. Moreover, if both
$P_\bullet:=M_\bullet\otimes_{V^a}A_\bullet$ and
$P_\bullet\otimes_{V^a}N_\bullet$ are almost essentially
constant with limits $P:=M\otimes_{V^a}A$ and respectively
$P\otimes_{V^a}N$, then
$M_\bullet\otimes_{V^a}A_\bullet\otimes_{V^a}N_\bullet$ will
be almost essentially constant with the sought limit. Hence,
we are further reduced to showing the proposition in case
$A_\bullet$ is the constant functor with $A_i=V^a$ for every
$i\in\Ob(I)$. Now, let $\tau_\bullet:c_M\Rightarrow M_\bullet$
be the universal cone. We have right exact sequences :
$$
K_i:=(\Ker\,\tau_i)\otimes_{V^a}N_i\to\Ker(\tau_i\otimes_{V^a}N_i)
\to T_i\to 0
\qquad
\text{for every $i\in\Ob(I)$}
$$
where $T_i$ is the cokernel of the natural morphism
$\Tor^{V^a}_1(M_i,N_i)\to\Tor^{V^a}_1(\Coker\,\tau_i,N_i)$.
Since, by assumption, both $\Ker\,\tau_\bullet$ and
$\Coker\,\tau_\bullet$ are almost essentially zero, it is
easily seen that the same holds for the induced functors
$K_\bullet$ and $T_\bullet$, and then the same follows for the
functor $\Ker(\tau_\bullet\otimes_{V^a}N_\bullet)$, by lemma
\ref{lem_a.ess.0-is-Serre}(i). Likewise,
$\Coker(\tau_\bullet\otimes_{V^a}N_\bullet)$ is the functor
$(\Coker\,\tau_\bullet)\otimes_{V^a}N_\bullet$, which is again
almost essentially zero. Thus, $M_\bullet\otimes_{V^a}N_\bullet$
is isomorphic to $M\otimes_{V^a}N_\bullet$ in the category
$\cL(I,(V,\fm)^a\Mod)$.
Lastly, let $F:(V,\fm)^a\Mod\to(V,\fm)^a\Mod$ be the functor
such that $FX:=M\otimes_{V^a}X$ for every $V^a$-module $X$;
since $N_\bullet$ is almost essentially constant, the same
holds for $FN_\bullet=M\otimes_{V^a}N_\bullet$, and its limit
is $M\otimes_{V^a}N$, by proposition
\ref{prop_apply-F-V-linear}(i,iii), whence the assertion.
\end{proof}

\begin{theorem}\label{th_alm-ess-cnst-inv-sys}
Let $(V,\fm)$ be a basic setup, $I$ a small cofiltered category,
and consider an almost essentially constant functor
$A_\bullet:I\to(V,\fm)^a\Alg$. Denote by $A$ the limit of $A_\bullet$,
and let $\tau^A_\bullet:c_A\Rightarrow A_\bullet$ be a universal cone;
let also $M$ be an $A$-module, $B$ an $A$-algebra, $r\in\N$ and
$$
A_\bullet\Mod:I\to\bCat
\quad
i\mapsto A_i\Mod
\qquad\text{and}\qquad
A_\bullet\Alg:I\to\bCat
\quad
i\mapsto A_i\Alg
$$
the pseudo-functors naturally induced by $A_\bullet$.
The following holds :

{\em (i)}\ \
The cone $\tau^A_\bullet$ induces equivalences of categories
$$
A\Alg\isom\Pslim{I}A_\bullet\Alg
\qquad
A\Mod\isom\Pslim{I}A_\bullet\Mod.
$$

{\em (ii)}\ \
The $A$-module $M$ is almost finitely generated (resp. almost
finitely presented, resp. flat, resp. faithfully flat resp.
almost projective) if and only if the same holds for the
$A_i$-module $A_i\otimes_AM$, for every $i\in\Ob(I)$.

{\em (iii)}\ \
Suppose that $(V,\fm)$ fulfills condition {\em($\bB$)} of\/
\cite[\S2.1.6]{Ga-Ra}. Then $M$ is almost projective of almost
finite rank (resp. of finite rank $\leq r$) if and only if the
same holds for the $A_i$-module $A_i\otimes_AM$, for every
$i\in\Ob(I)$.

{\em (iv)}\ \
The $A$-algebra $B$ is weakly unramified (resp. unramified,
resp. weakly \'etale, resp. \'etale) if and only if the same
holds for the $A_i$-algebra $A_i\otimes_AB$, for every $i\in\Ob(I)$.
\end{theorem}
\begin{proof}(i): Let $\pi:(V,\fm)^a\tdu\AlgMod\to(V,\fm)^a\Alg$
be the natural projection functor. Then $\Pslim{I}A_\bullet\Mod$
is the subcategory of $\bFun(I,(V,\fm)^a\tdu\AlgMod)$ whose
objects are the functors
$$
(A_\bullet,M_\bullet):I\to(V,\fm)^a\tdu\AlgMod
\qquad
i\mapsto(A_i,M_i)
$$
with $\pi\circ(A_\bullet,M_\bullet)=A_\bullet$, and such that, with
the notation of \eqref{prop_inverse-syst-of-AlgMods}, the
morphism $g_\phi:A_i\otimes_{A_j}M_j\to M_i$ is an isomorphim of
$A_i$-modules, for every morphism $\phi:j\to i$ in $I$. The
morphisms $\beta_\bullet:(A_\bullet,M_\bullet)\to(A_\bullet,M'_\bullet)$
are the natural transformations such that
$\pi*\beta_\bullet=\one_{A_\bullet}$.

Under this identification, the functor
$(A_\bullet,A_\bullet\otimes_A-):A\Mod\to\Pslim{I}A_\bullet\Mod$
of (i) assigns to every $A$-module $M$ the functor
$$
(A_\bullet,A_\bullet\otimes_AM):I\to(V,\fm)^a\tdu\AlgMod
\qquad
i\mapsto(A_i,A_i\otimes_AM).
$$
Conversely, following \eqref{prop_inverse-syst-of-AlgMods},
with every functor $(A_\bullet,M_\bullet)$ as in the foregoing, we
associate a functor $M_\bullet:I\to A\Mod$, and let $M$ be the limit
of $M_\bullet$. From proposition \ref{prop_return-of-AlMod}(iii)
we obtain a natural isomorphism
$$
(A_\bullet,A_\bullet\otimes_AM)\isom(A_\bullet,M_\bullet).
$$
This shows that $(A_\bullet,A_\bullet\otimes_A-)$ is
essentially surjective. Next, let $M,N$ be two $A$-modules;
according to proposition \ref{prop_apply-F-V-linear}(i),
the kernel and cokernel of
$\tau^A_\bullet\otimes_AN:c_N\Rightarrow A_\bullet\otimes_AN$ are
almost essentially zero, and therefore $\tau^A_\bullet\otimes_AN$
is a universal cone. There follow natural isomorphisms
$$
\Hom_A(M,N)\isom\lim_I\Hom_A(M,A_\bullet\otimes_AN)\isom
\lim_I\Hom_{A_i}(A_\bullet\otimes_AM,A_\bullet\otimes_AN)
$$
and notice that the latter limit is naturally identified
with the set of morphisms
$(A_\bullet,A_\bullet\otimes_AM)\to(A_\bullet,A_\bullet\otimes_AN)$
in $\Pslim{I}A_\bullet\Mod$. Thus $(A_\bullet,A_\bullet\otimes_A-)$
is an equivalence.

(ii): Suppose first that $M_i:=A_i\otimes_AM$ is almost finitely
generated for every $i\in\Ob(I)$. Let $N_\bullet:\Lambda\to A\Mod$
be any functor from a small filtered category $\Lambda$; denote
by $N$ the colimit of $N_\bullet$, and let
$\tau^N_\bullet:N_\bullet\Rightarrow c_N$ be a universal cocone.
We have two $V$-linear functors
$$
\begin{aligned}
F&\,:A\Mod\to A\Mod & \qquad &
P\mapsto\colim_\Lambda\,\Alhom_A(M,P\otimes_AN_\bullet) \\
G&\,:A\Mod\to A\Mod & \qquad &
P\mapsto\Alhom_A(M,P\otimes_AN)
\end{aligned}
$$
and $\tau^N_\bullet$ induces a natural transformation
$$
\beta_\bullet:F\Rightarrow G
\qquad
P\mapsto\colim_\Lambda\,\Alhom_A(M,P\otimes_A\tau^N_\bullet):FP\to GP
$$
and according to \cite[Prop.2.3.16(i)]{Ga-Ra}, it suffices to
check that $\beta_A$ is a monomorphism. However, we have a
commutative diagram of functors :
$$
\xymatrix@C+20pt{
c_{FA} \ar[r]^-{c_{\beta_A}} \ar[d]_{F*\tau^A_\bullet} &
c_{GA} \ar[d]^{G*\tau^A_\bullet} \\
F\circ A_\bullet \ar[r]^-{\beta_\bullet*A_\bullet} & G\circ A_\bullet
}$$
and by propositions \ref{prop_a-Mittag-Leffler}(ii) and
\ref{prop_apply-F-V-linear}(i), the vertical arrows of
the diagram are both isomorphisms in the category
$\cL(I,A\Mod)$. We claim that $\beta_{A_i}:FA_i\to GA_i$
is a monomorphism for every $i\in\Ob(I)$. Indeed :
$$
FA_i=\colim_\Lambda\,\Alhom_{A_i}(M_i,A_i\otimes_AN_\bullet)
\qquad\text{and}\qquad
GA_i=\Alhom_{A_i}(M_i,A_i\otimes_AN)
$$
and since by assumption the $A_i$-module $M_i$ is almost
finitely generated, the assertion follows from
\cite[Prop.2.3.16(i)]{Ga-Ra}. We deduce that $c_{\beta_A}$
is also a monomorphism in $\cL(I,A\Mod)$, {\em i.e.}
$\beta_A$ is a monomorphism, as required.

One argues likewise, in case $M_i$ is almost finitely
presented for every $i\in\Ob(I)$, invoking
\cite[Prop.2.3.16(ii)]{Ga-Ra} to see that in this case
$\beta_{A_i}$ is an isomorphism for every such $i$, whence
$\beta_A$  is an isomorphism, so that $M$ is almost finitely
presented, again by \cite[Prop.2.3.16(ii)]{Ga-Ra}.

Next, suppose that $M_i$ is an almost projective $A_i$-module
for every $i\in\Ob(I)$, and let $f:N'\to N$ be an epimorphism
of $A$-modules. We attach to the pair $(M,N)$ the $V$-linear
functor
$$
H_{M,N}:A\Mod\to A\Mod
\qquad
P\mapsto\Alhom_A(P\otimes_AM,P\otimes_AN)
$$
and define likewise $H_{M,N'}:A\Mod\to A\Mod$. The morphism
$f$ induces a natural transformation
$H_{M,f}:H_{M,N}\Rrightarrow H_{M,N'}$, whence a commutative
diagram of functors :
$$
\xymatrix@C+30pt{
c_{H_{M,N}(A)} \ar[r]^-{c_{H_{M,f}(A)}} \ar[d]_{H_{M,N}*\tau^A_\bullet} &
c_{H_{M,N'}(A)} \ar[d]^{H_{M,N'}*\tau^A_\bullet} \\
H_{M,N}\circ A_\bullet \ar[r]^-{H_{M,f}*A_\bullet} &
H_{M,N'}\circ A_\bullet
}$$
whose vertical arrows are isomorphisms in $\cL(I,A\Mod)$,
by virtue of propositions \ref{prop_a-Mittag-Leffler}(ii) and
\ref{prop_apply-F-V-linear}(i); moreover, $H_{M,N'}*\tau^A_i$
is an epimorphism for every $i\in\Ob(I)$, since $M_i$ is an
almost projective $A_i$-module. As in the previous case,
we deduce that $H_{M,f}(A)$ is an epimorphism of $A$-modules,
so $M$ is almost projective.

Suppose next that $M_i$ is flat for every $i\in\Ob(I)$, and
let $f:N\to N'$ be a monomorphism of $A$-modules. Set
$N_i:=A_i\otimes_AN$ and $N'_i:=A_i\otimes_AN'$ for every
$i\in\Ob(I)$, and let $\bar N_i:=\Img(A_i\otimes_Af)\subset N'_i$.
Clearly the rule : $i\mapsto\bar N_i$ defines a subfunctor
$\bar N_\bullet$ of $A_\bullet\otimes_AN'$. By proposition
\ref{prop_apply-F-V-linear}(i), the kernel of
$A_\bullet\otimes_Af:A_\bullet\otimes_AN\Rightarrow A_\bullet\otimes_AN'$
is almost essentially zero, hence the induced natural transformation
$\pi_\bullet:A_\bullet\otimes_AN\to\bar N_\bullet$ is an isomorphism in
$\cL(I,A\Mod)$. By the same token, the same holds for the natural
transformation
$\tau^A_\bullet\otimes_AN:c_N\Rightarrow A_\bullet\otimes_AN$,
and therefore also for the composition
$\mu_\bullet:=\pi_\bullet\odot(\tau^A_\bullet\otimes_AN):
c_N\Rightarrow\bar N_\bullet$. We obtain therefore a commutative
diagram of functors :
$$
\xymatrix@C+20pt{
c_{M\otimes_AN} \ar[r]^-{c_{M\otimes_Af}} \ar[d]_{M\otimes_A\mu_\bullet} &
c_{M\otimes_AN'} \ar[d]^{M\otimes_A\tau^A_\bullet\otimes_AN} \\
M\otimes_A\bar N_\bullet \ar[r]^-{} & M\otimes_AA_i\otimes_AN'
}$$
whose vertical arrows are, as usual, isomorphisms in
$\cL(I,A\Mod)$. Let $j_i:\bar N_i\to N'_i$ be
the inclusion; we have natural isomorphisms
$M\otimes_A\bar N_i\isom M_i\otimes_{A_i}\bar N_i$ and
$M\otimes_AA_i\otimes_AN'\isom M_i\otimes_{A_i}N'_i$ for
every $i\in\Ob(I)$, that identify the bottom horizontal
arrow with the morphism $M\otimes_{A_i}j_i$. The latter
is a monomorphism for every $i\in\Ob(I)$, since $M_i$
is a flat $A_i$-module; as usual it follows that
$M\otimes_Af$ is a monomorphism as well, so $M$ is flat.

Lastly, suppose that $M_i$ is a faithfully flat $A_i$-module
for every $i\in\Ob(I)$; by the foregoing, we know already
that $M$ is a flat $A$-module. Then, let $X$ be any $A$-module
such that $M\otimes_AX=0$; it follows that
$M_i\otimes_{A_i}(A_i\otimes_AX)=0$ for every $i\in\Ob(I)$,
whence $A_i\otimes_AX=0$ for every such $i$. But by the
foregoing, $X$ is the limit of the functor $A_\bullet\otimes_AX$,
whence $X=0$; this shows that $M$ is faithfully flat.

(iii): Suppose that $M_i$ is an almost projective
$A_i$-module of almost finite rank, for every $i\in\Ob(I)$.
By the foregoing, we know already that $M$ is an almost
projective $A$-module, and it remains to check that it
is of almost finite rank. To this aim, define the category
$\cI$ and the functors $\pi:\cI\to\sMorph(I)$ and
$\bar A_\bullet:\sMorph(I)\to A\Alg$ as in the proof of
proposition \ref{prop_return-of-AlMod}; recall that
$\pi$ is final, and $\bar A_\bullet\circ\pi$ is a Cauchy
functor (claims \ref{cl_pi-and-delta-final} and
\ref{cl_Cauchification}). Moreover, the cone
$\bar\tau{}^A_\bullet:c_A\Rightarrow\bar A_\bullet$ deduced
from $\tau^A_\bullet$, is universal, so the same holds
for $\bar\tau{}^A_\bullet*\pi$. It follows that
$\Ker\,(\bar\tau{}^A_\bullet*\pi)$ and
$\Coker\,(\bar\tau{}^A_\bullet*\pi)$ are null functors
(remark \ref{rem_Cauchy-functors}(ii)). Now, let
$\fm_0\subset\fm$ be a subideal of finite type; by the
foregoing there exists $(\phi:j\to i,\fm_1)\in\Ob(\cI)$
such that the kernel and cokernel of
$\bar\tau{}^A_\phi:A\to\bar A_\phi:=\Img\,A_\phi$ are both
annihilated by $\fm_0$. On the other hand, since $M_j$
has almost finite rank, there exists $n\in\N$ such that
$\fm_0\cdot\Lambda^n_{A_j}M_j=0$, and then we have as well
$\fm_0\cdot(\bar A_\phi\otimes_{A_j}\Lambda^n_{A_j}M_j)=0$. But
by remark \ref{rem_V-lin-abelian-cat}(vi), the kernel of
$$
\bar\tau{}^A_\phi\otimes_A\Lambda^n_AM:
\Lambda^n_AM\to\bar A_\phi\otimes_A\Lambda^n_AM\isom
\bar A_\phi\otimes_{A_j}\Lambda^n_{A_j}M_j
$$
is annihilated by $\fm_0^2$. We conclude that
$\fm_0^3\cdot\Lambda^n_AM=0$, whence the contention.
Lastly, if $M_i$ has rank $\leq r$ for every $i\in\Ob(I)$,
the same argument shows that $\fm\cdot\fm_0^2\cdot\Lambda^r_AM=0$
for every such $\fm_0$, {\em i.e.} $\Lambda^r_AM=0$, so
$M$ has rank $\leq r$.

(iv): Set $C:=B\otimes_AB$, and let $\mu:C\to B$ be the
multiplication law of the $A$-algebra $B$; the functor
$C_\bullet:=C\otimes_AA_\bullet:I\to(V,\fm)^a\Alg$ is still almost
essentially constant, and
$C\otimes_A\tau^A_\bullet:c_C\Rightarrow C_\bullet$ is an isomorphism
in $\cL(I,A\Mod)$ (proposition \ref{prop_apply-F-V-linear}(i,iii)).
Suppose now that $B_i:=B\otimes_AA_i$ is a weakly unramified
$A_i$-algebra, {\em i.e.} that the morphism
$\mu\otimes_AA_i:C_i\to B_i$ is flat for every $i\in\Ob(I)$.
By (ii), it follows that $\mu$ is flat, {\em i.e.} $B$ is a
weakly unramified $A$-algebra. One argues likewise in case
$B_i$ is an unramified (resp. weakly \'etale, resp. \'etale)
$A_i$-algebra for every $i\in\Ob(I)$.
\end{proof}

\subsection{Quasi-coherent sheaves of almost modules and
almost rings}\label{subsec_qcoh-almost-rings}
Throughout this section, we fix a basic setup $(V,\fm)$ (see
\cite[\S2.1.1]{Ga-Ra}) and we set $S:=\Spec\,V$; for every
$S$-scheme $X$ we may consider the sheaf $\cO_{\!X}^a$ of
almost algebras on $X$, and we refer to \cite[\S5.5]{Ga-Ra}
for the definition of quasi-coherent $\cO^a_{\!X}$-modules
and algebras. However, whereas in \cite{Ga-Ra} it was assumed
that $\tilde\fm:=\fm\otimes_V\fm$ is a flat $V$-module,
in this section the basic setup can be arbitrary, except
where it is explicitly said otherwise. This of course
requires some care when quoting from \cite{Ga-Ra}; however,
it turns out that -- thanks to the work done in section
\ref{sec_non-flat} -- most of the results in \cite{Ga-Ra} do
extend {\em verbatim} to the case of a general setup. The main
exception is the theory of the finite and almost finite rank
of almost projective modules, that relies on the existence
of a well behaved exterior power functor, which is available
only if the basic setup satisfies some minimal conditions
(see \eqref{subsec_drop_B}). In any case, whenever we need to
import some theorem from \cite{Ga-Ra}, we shall comment on its
range of validity.

\begin{definition}\label{def_complement-alqcoh}
Let $X$ be an $S$-scheme, $\cA$ a quasi-coherent
$\cO_{\!X}^a$-algebra, and $\cF$ an $\cA$-module which
is quasi-coherent as an $\cO^a_{\!X}$-module.

(i)\ \
$\cF$ is said to be an $\cA$-module {\em of almost finite type}
(resp. {\em of almost finite presentation}, resp {\em flat},
resp. {\em faithfully flat}) if, for every affine open subset
$U\subset X$, the $\cA(U)$-module $\cF(U)$ is almost finitely
generated (resp. almost finitely presented, resp. flat, resp.
faithfully flat).

(ii)\ \
$\cF$ is said to be an {\em almost coherent\/} $\cA$-module
if it is an $\cA$-module of almost finite type, and for every
open subset $U\subset X$, every quasi-coherent $\cA_{|U}$-submodule
of $\cF_{|U}$ of almost finite type, is almost finitely presented.

(iii)\ \
We say that $\cF$ is a {\em torsion-free\/} $\cO^a_{\!X}$-module
if we have $\Ker\,b\cdot\one_{\cF(U)}=0$, for every affine open
subset $U\subset X$ and every regular element $b\in\cO_{\!X}(U)$.

(iv)\ \
Suppose that $\fm$ satisfies condition $(\bB)$ of
\cite[\S2.1.6]{Ga-Ra}. Then we say that $\cF$ is an $\cA$-module
{\em of almost finite rank} (resp. {\em of finite rank}) if, for
every affine open subset $U\subset X$, the $\cA(U)$-module
$\cF(U)$ is almost finitely generated projective of almost
finite (resp. finite) rank.
\end{definition}

\begin{remark}\label{rem_global-trace-morph}
In the situation of definition \ref{def_complement-alqcoh} :

(i)\ \
We denote by
$$
\cEnd_\cA(\cF)
$$
the $\cO_{\!X}$-module of {\em $\cA$-linear endomorphisms of $\cF$},
defined by the rule : $U\mapsto\End_{\cA(U)}(\cF(U))$ for every open
subset $U\subset X$. It is easily seen that $\cEnd_\cA(\cF)^a$ is an
$\cA$-module.

(ii)\ \
Suppose that $\cF$ is a flat and almost finitely presented
$\cA$-module. Then there exists a {\em trace morphism}
$$
\tr_{\cF/\cA}:\cEnd_\cA(\cF)\to\cA
$$
that, on every affine open subset $U\subset X$, induces the
trace morphism $\tr_{\cF(U)\cA(U)}$ of the almost finitely
generated projective $\cA(U)$-module $\cF(U)$ (details left
to the reader : see \cite[\S4.1]{Ga-Ra}, which does not
depend on any assumption on the basic setup).
\end{remark}

\sset\subsubsection{}\label{subsec_almost-up-down}
Let $f:Y\to X$ be any morphism of $V$-schemes. The usual
functor $f^*$ for quasi-coherent $\cO_{\!X}$-module
admits a variant for quasi-coherent $\cO_{\!X}^a$-modules;
namely, we define
$$
f^*:\cO_{\!X}^a\Mod_\qcoh\to\cO^a_Y\Mod_\qcoh
\qquad
\cF\mapsto(f^*\cF_!)^a.
$$
If $f$ is also quasi-compact and quasi-separated, we
have as well a direct image functor
$$
f_*:\cO_Y^a\Mod_\qcoh\to\cO_{\!X}^a\Mod_\qcoh
\qquad
\cG\mapsto(f_*\cG_!)^a
$$
(\cite[Ch.I, Prop.9.2.1]{EGAI}). Moreover, if
$\cF$ (resp. $\cG$) is any quasi-coherent
$\cO_{\!X}$-module (resp. $\cO_Y$-module), it is easily
seen that the natural morphism
$$
f^*(\cF^a)\to(f^*\cF)^a
\qquad
\text{resp.\ $f_*(\cG^a)\to(f_*\cG)^a$}
$$
is an isomorphism. It follows that, if $Z\subset X$
is any constructible closed subset, there is a
well defined functor of sections with support in $Z$ :
$$
\underline\Gamma_Z:\cO_{\!X}\Mod_\qcoh\to\cO_{\!X}\Mod_\qcoh
\qquad
\cF\mapsto(\underline\Gamma_Z\cF_!)^a
$$
(see \eqref{subsec_case-of-schs}) and again, for every
quasi-coherent $\cO_{\!X}$-module $\cF$, the natural map
$$
\underline\Gamma_Z(\cF^a)\to(\underline\Gamma_Z\cF)^a
$$
is an isomorphism (details left to the reader).

\begin{lemma}\label{lem_reduced-or-integral}
Suppose that $\fm$ fulfills condition $(\bB)$, and let $X$
be a reduced $S$-scheme, $\cF$ a flat quasi-coherent
$\cO^a_X$-module of almost finite type. We have :
\begin{enumerate}
\item
If $X$ is reduced and $\cF$ is an $\cO^a_{\!X}$-module of almost
finite presentation, then $\cF$ is an $\cO^a_{\!X}$-module of
almost finite rank.
\item
If $X$ is integral, then $\cF$ is an $\cO^a_{\!X}$-module
of finite rank.
\end{enumerate}
\end{lemma}
\begin{proof}(i): We may assume that $X$ is affine, say
$X=\Spec\,R$ for a reduced $V$-algebra $R$, and $\cF=P^\sim$
for a flat and almost finitely presented $R^a$-module $P$.
Then $P$ is almost projective (\cite[Prop.2.4.18(ii)]{Ga-Ra}),
and we need to check that it has almost finite rank.
Thus, let $\eps\in\fm$ be any element and set
$R':=R[\eps^{-1}]$; then $P':=R'\otimes_{R^a}P$ is a projective
$R'$-module of finite rank, so there exists $n\in\N$ such
that $R'\otimes_R\Lambda^n_RP=\Lambda^n_{R'}P'=0$. Since
$\Lambda^n_RP$ is also almost finitely presented (see
\cite[\S4.3.1]{Ga-Ra}), we deduce easily that there exists
$N\in\N$ such that $\eps^N\cdot\Lambda^n_RP=0$, and since
moreover $\Lambda^n_RP$ is a flat $R^a$-module and $R$ is
reduced, we have
$$
\Ann_{\Lambda^n_RP}(\eps^N)=\Ann_R(\eps^N)\cdot\Lambda^n_RP=
\Ann_R(\eps)\cdot\Lambda^n_RP
$$
{\em i.e.} $\eps\cdot\Lambda^n_RP=0$, whence the contention.

(ii): Let $\eta$ be the generic point of $X$; if
$\kappa(\eta)^a=0$, we have $\cO^a_{\!X}=0$ as well, and there
is nothing to show. If $\kappa(\eta)^a\neq 0$, the almost
structure on $\kappa(\eta)$ is the trivial one (the ``classical
limit'' case of \cite[Ex.2.1.2(ii)]{Ga-Ra}); in this case,
clearly $\cF_\eta$ is a free $\kappa(\eta)^a$-module of finite
rank, and let $r$ be this rank.
From \cite[Prop.2.4.19]{Ga-Ra}, it follows easily that $\cF$
is almost finitely presented. Moreover, the exterior powers of
$\cF$ are still flat $\cO^a_{\!X}$-modules; then the $r$-th
exterior power vanishes, since it vanishes at the generic
point of $X$.
\end{proof}

\begin{lemma}\label{lem_complem-qcoh}
Let $f:Y\to X$ a faithfully flat and quasi-compact morphism
of $S$-schemes, $\cA$ a quasi-coherent $\cO^a_{\!X}$-algebra,
$\cF$ a quasi-coherent $\cA$-module, and $r\in\N$. Then :
\begin{enumerate}
\item
The $\cA$-module $\cF$ is of almost finite type (resp. of
almost finite presentation, resp. flat, resp. faithfully flat)
if and only if the same holds for the $f^*\!\cA$-module $f^*\!\cF$.
\item
Suppose that $\fm$ fulfills condition $(\bB)$. Then the
$\cA$-module $\cF$ is of almost finite rank (resp. of finite
rank $\leq r$) if and only if the same holds for the
$f^*\!\cA$-module $f^*\!\cF$.
\item
If the $f^*\!\cA$-module $f^*\!\cF$ is almost coherent, then the
same holds for the $\cA$-module $\cF$.
\end{enumerate}
\end{lemma}
\begin{proof} The assertions are local on $X$, so
we may assume that $X$ is affine; then $Y$ admits
a finite covering $(U_i~|~i\in I)$ consisting of
affine open subsets, and we may further reduce to
the case where $Y$ is the disjoint union of the
schemes $U_i$, so $Y$ is affine as well.

(i): For the conditions ``almost finite type'', ``almost finite
presentation'', and ``flat'', the assertion follows from
\cite[Rem.3.2.26(ii)]{Ga-Ra}, which holds for any basic setup.
It remains to check that if $A\to B$ is a faithfully flat map
of $V$-algebras, and $M$ is any $A$-module, such that
$(B\otimes_AM)^a$ is a faithfully flat $B^a$-module, then $M^a$
is a faithfully flat $A^a$-module. To this aim, let $X$ be any
$A$-module such that $(M\otimes_AX)^a=0$; then
$(B\otimes_AX)^a\otimes_{B^a}(B\otimes_AX)^a=0$, and therefore
$(B\otimes_AX)^a=0$, {\em i.e.} $\tilde\fm\otimes_VB\otimes_AX=0$,
whence $\tilde\fm\otimes_VX=0$, {\em i.e.} $X^a=0$. Since we know
already that $M^a$ is a flat $A^a$-module, the contention follows.

(ii): It suffices to apply \cite[Rem.3.2.26(iii)]{Ga-Ra},
and recall that exterior powers commute with any base changes :
the details shall be left to the reader.

(iii) follows easily from (i).
\end{proof}

\begin{lemma}\label{lem_al-coherence}
Let $X$ be any $S$-scheme.
\begin{enumerate}
\item
The full subcategory $\cO^a_{\!X}\Mod_\mathrm{acoh}$ of
$\cO^a_{\!X}\Mod_\qcoh$ consisting of all almost coherent
modules is abelian and closed under extensions. (More
precisely, the embedding
$\cO^a_{\!X}\Mod_\mathrm{acoh}\to\cO^a_{\!X}\Mod_\qcoh$
is an exact functor.)
\item
If $X$ is a coherent scheme, then every quasi-coherent
$\cO^a_{\!X}$-module of almost finite presentation is almost
coherent.
\end{enumerate}
\end{lemma}
\begin{proof} Both assertions are local on $X$, hence we may
assume that $X$ is affine, say $X=\Spec\,R$.

(i): Let $f:M\to N$ be a morphism of almost coherent $R^a$-modules.
It is clear that $\Coker\,f$ is again almost coherent. Moreover,
$f(M)\subset N$ is almost finitely generated, hence almost finitely
presented; then by \cite[Lemma 2.3.18]{Ga-Ra}$, \Ker\,f$ is
almost finitely generated, and so it is almost coherent as well.
Similarly, using \cite[Lemma 2.3.18]{Ga-Ra} one sees that
$\cO^a_{\!X}\Mod_\mathrm{acoh}$ is closed under extensions.

(ii): Let $M$ be an almost finitely presented $R^a$-module;
according to \cite[Cor.2.3.13]{Ga-Ra}, for every finitely
generated ideal $\fm_0\subset\fm$ there exist a finitely
presented $R$-module $N$ and an $R$-linear map $\phi:N\to M_*$
whose kernel and cokernel are annihilated by $\fm_0$. Likewise,
if $M'\subset M$ is an almost finitely generated
$R^a$-submodule, we may find finitely many almost
elements $x_1,\dots,x_n\in M'_*$ that generate an
$R$-module containing $\fm_0 M'_*$. Let $a_1,\dots,a_k$
be a finite system of generators for $\fm_0$; for each
$i\leq k$ and $j\leq n$ we may find $y_{ij}\in N$ such
that $\phi(y_{ij})=a_ix_j$. Denote by $N'$ the $R$-submodule
of $N$ generated by $(y_{ij}~|~i\leq k,\ j\leq n)$. Then
$\phi$ restricts to an $R$-linear map $N'\to M'_*$ whose
kernel is annihilated by $\fm_0$ and whose cokernel is
annihilated by $\fm_0^2$. Since $R$ is coherent, $N'$
is finitely presented, and since $\fm_0$ is arbitrary,
we conclude that $M'$ is almost finitely presented,
as stated.
\end{proof}

\begin{lemma}\label{lem_alm-faithful-flatness}
Let $X$ be an $S$-scheme, $\phi:\cA\to\cB$ a morphism of
quasi-coherent $\cO_X^a$-algebras, and
$\cB_\bullet:=(\cB_\lambda~|~\lambda\in\Lambda)$ a filtered
system of quasi-coherent and faithfully flat $\cA$-algebras.
We have :

{\em(i)}\ \
$\phi$ is faithfully flat if and only if $\phi$ is injective
and $\Coker\,\phi$ is a flat $\cA$-module.

{\em(ii)}\ \
The colimit of $\cB_\bullet$ is a faithfully flat $\cA$-algebra. 
\end{lemma}
\begin{proof} It is easily seen that (i)$\Rightarrow$(ii).
In order to show (i), we may assume that $X$ is affine,
in which case $\phi$ is the morphism of quasi-coherent
$\cO_X^a$-algebras arising from a morphism $f:A\to B$
of $\cO_X(X)^a$-algebras. Now, let $M$ be any $A$-module;
if $f$ is injective, the short exact sequence of $A$-modules
$0\to A\to B\to\Coker\,f\to 0$ induces a long exact sequence
$$
0\to\Tor^A_1(B,M)\to\Tor^A_1(\Coker\,f,M)\to
M\xrightarrow{f\otimes_AM}B\otimes_AM\to
\Coker\,f\otimes_AM\to 0.
$$
Then, if $\Coker\,f$ is a flat $A$-module, we deduce that
$f\otimes_AM$ is injective and $\Tor^A_1(B,M)=0$ for every
$A$-module $M$, so $B$ is a faithfully flat $A$-algebra.
Conversely, if the latter holds, then $f$ and $f\otimes_MA$
are both injective, hence the foregoing short exact sequence
shows that $\Tor^A_1(\Coker\,f,M)=0$ for every $A$-module $M$,
so that $\Coker\,f$ is a flat $A$-module.
\end{proof}

\begin{definition}\label{def_complement-alqcoh-alg}
Let $X$ be an $S$-scheme, and $\cA$ a quasi-coherent
$\cO^a_{\!X}$-algebra.

(i)\ \
We say that a quasi-coherent $\cA$-algebra $\cB$ is
{\em almost finite\/} (resp. {\em weakly unramified},
resp. {\em weakly \'etale}, resp. {\em unramified}, resp.
{\em \'etale}) if, for every affine open subset $U\subset X$,
the $\cA(U)$-algebra $\cB(U)$ is almost finite (resp. weakly
unramified, resp. weakly \'etale, resp. unramified, resp. \'etale) .

(ii)\ \
We define the $\cO_{\!X}$-algebra $\cA_{!!}$ by the short
exact sequence :
$$
\tilde\fm\otimes_V\cO_{\!X}\to\cO_{\!X}\oplus\cA_!\to\cA_{!!}\to 0
$$
analogous to \cite[\S2.2.25]{Ga-Ra}. This agrees with the definition
of \cite[\S3.3.2]{Ga-Ra}, corresponding to the  basic setup
$(\cO_{\!X},\fm\cO_{\!X})$ relative to the Zariski topos of $X$;
in general, this is {\em not\/} the same as forming the algebra
$\cA_{!!}$ relative to the basic setup $(V,\fm)$. The reason why
we prefer the foregoing version of $\cA_{!!}$, is explained by
the following lemma \ref{lem_about-integral-closure}.
\end{definition}

\begin{lemma}\label{lem_about-integral-closure}
Let $X$ be an $S$-scheme, $\cA$ a quasi-coherent
$\cO^a_{\!X}$-algebra. Then $\cA_{!!}$ is a quasi-coherent
$\cO_{\!X}$-algebra, and if $\cA$ is almost finite, then
$\cA_{!!}$ is an integral $\cO_{\!X}$-algebra.
\end{lemma}
\begin{proof} First, $\cA_!$ is a quasi-coherent $\cO_{\!X}$-module
(\cite[\S5.5.4]{Ga-Ra}), hence the same holds for $\cA_{!!}$. It
remains to show that if $\cA$ is almost finite, then
$\cA_{!!}(U)=\cA(U)_{!!}$ is an integral $\cO_{\!X}(U)$-algebra
for every affine open subset $U\subset X$.
Since $\cA(U)_{!!}=\cO_{\!X}(U)+\fm\cA(U)_{!!}$, we need only
show that every element of $\fm\cA(U)_{!!}$ is integral over
$\cO_{\!X}(U)$. However, by adjunction we get a map :
$$
\cA(U)_{!!}\to A_U:=\cO^a_{\!X}(U)_*+\fm\cA(U)_*
$$
whose kernel is annihilated by $\fm$, and according to
\cite[Lemma 5.1.13(i)]{Ga-Ra}, $A_U$ is integral over
$\cO^a_{\!X}(U)_*$. Let $a\in\fm\cA(U)_{!!}$ be any element,
and $\bar a$ its image in $\fm\cA(U)_*$; we can then find
$b_0,\dots,b_n\in\cO_{\!X}(U)_*$ such that
$\bar a^{n+1}+\sum^n_{i=0}b_i\bar a{}^i=0$ in $A_U$. It follows
that $(\eps a)^{n+1}+\sum^n_{i=0}\eps^{n+1-i}b_i(\eps a)^i=0$
in $\cA(U)_{!!}$, for every $\eps\in\fm$. Since $\eps^{n+1-i}b_i$
lies in the image of $\cO_{\!X}(U)$ for every $i\leq n$, the
claim follows easily.
\end{proof}

\begin{remark}\label{rem_idemp-and-traces}
With the notation of definition \ref{def_complement-alqcoh-alg}
we have :

(i)\ \
$\cB$ is an unramified $\cA$-algebra if and only if, for every
affine open subset $U\subset X$, there exists an idempotent
element $e_{\cB(U)/\cA(U)}\in(\cB\otimes_\cA\cB)(U)_*$, uniquely
characterized by the following conditions :
\begin{enumerate}
\alphaenu
\item
$\mu_{\cB/\cA}(e_{\cB(U)/\cA(U)})=1$, where
$\mu_{\cB/\cA}:\cB\otimes_{\cA}\cB\to\cB$ is the multiplication
morphism of the $\cA$-algebra $\cB$.
\item
$e_{\cB(U)/\cA(U)}\cdot\Ker\,\mu_{\cB/\cA}(U)=0$
\end{enumerate}
(\cite[Lemma 3.1.4]{Ga-Ra}). It is easily seen that, on $U'\subset U$,
the element $e_{\cB(U)/\cA(U)}$ restricts to $e_{\cB(U')/\cA(U')}$;
hence, the system $(e_{\cB(U)/\cA(U)}~|~U\subset X)$, for $U$ ranging
over all the affine open subsets of $X$, determines a global section
$$
e_{\cB/\cA}\in\Gamma(X,\cB\otimes_\cA\cB)_*
$$
which we call the {\em diagonal idempotent\/} of the $\cA$-algebra
$\cB$ (see \cite[\S5.5.4]{Ga-Ra}).

(ii)\ \
In the situation of (i), let $\cA\to\cA'$ be any morphism
of quasi-coherent $\cO^a_{\!X}$-algebras, and set
$\cB':=\cA'\otimes_\cA\cB$. Then the induced morphism
$\cA'\to\cB'$ is unramified, and in view of the characterization
provided by (i), it is easily seen that $e_{\cB'/\cA'}$ is
the image in $\cB'\otimes_{\cA'}\cB'$ of $e_{\cB/\cA}$.

(iii)\ \
If $\cB$ is a flat and  almost finitely presented $\cA$-algebra,
there exists a {\em trace form}
$$
t_{\cB/\cA}:\cB\otimes_\cA\cB\to\cA
$$
that, on each open affine subset $U\subset X$, induces the
trace form $t_{\cB(U)/\cA(U)}$ of \cite[\S4.1.12]{Ga-Ra}
(details left to the reader).

(iv)\ \
Suppose that $\cB$ is a flat, unramified and almost finitely
presented $\cA$-algebra. Then the diagonal idempotent and the
trace form of $\cB$ are related by the identity
$$
t_{\cB/\cA}(e_{\cB/\cA}\cdot(1\otimes b))=b=
t_{\cB/\cA}(e_{\cB/\cA}\cdot(b\otimes 1))
$$
for every affine open subset $U\subset X$ and every
$b\in\cB(U)_*$ (\cite[Rem.4.1.17]{Ga-Ra}).
\end{remark}

\begin{definition}\label{def_alm-int-closure}
(i)\ \
Let $X$ be an $S$-scheme. For a monomorphism $\cR\subset\cS$
of quasi-coherent $\cO_{\!X}$-algebras, the integral closure
$\mathrm{i.c.}(\cR,\cS)$ of $\cR$ in $\cS$ is a quasi-coherent
$\cO_{\!X}$-algebra. For a monomorphism $\cA\to\cB$ of
$\cO_X^a$-algebras, the {\em integral closure of $\cA$ in $\cB$}
is
$$
\mathrm{i.c.}(\cA,\cB):=\mathrm{i.c.}(\cA_{!!},\cB_{!!})^a
$$
(see also \cite[Lemma 8.2.28]{Ga-Ra}). In view of lemma
\ref{lem_about-integral-closure}, this is a well defined
quasi-coherent $\cO^a_{\!X}$-algebra. It is characterized
as the unique $\cA$-subalgebra of $\cB$ such that :
$$
\mathrm{i.c.}(\cA,\cB)(U)=\mathrm{i.c.}(\cA(U),\cB(U))
$$
for every affine open subset $U\subset X$ (notation of
\cite[Def.8.2.27]{Ga-Ra}). We say that $\cA$ is {\em
integrally closed in $\cB$}, if $\cA=\mathrm{i.c.}(\cA,\cB)$.
We say that $\cB$ is an {\em integral\/} $\cA$-algebra
if $\mathrm{i.c.}(\cA,\cB)=\cB$.

(ii)\ \
Likewise, the weak normalization $\cC$ and the $p$-integral
closure $\cD$ of the image of $\cA_{!!}$ in $\cB_{!!}$ are
quasi-coherent $\cO_X$-algebras (see
\eqref{subsec_global-weak-norm}); then we define the
{\em weak normalization} (resp. the {\em $p$-integral
closure}) of $\cA$ in $\cB$ as the quasi-coherent
$\cO_X^a$-algebra $\cC^a$ (resp. $\cD^a$).
\end{definition}

\begin{remark}\label{rem_alm-int-closure}
(i)\ \
In the notation of definition \ref{def_alm-int-closure},
if $\cR'$ is the integral closure of $\cR$ in $\cS$, then
$\cR'^a$ is the integral closure of $\cR^a$ in $\cS^a$.
Indeed, we may assume that $X$ is affine, in which case the
assertion follows from \cite[Lemma 8.2.28]{Ga-Ra}.

(ii)\ \
It follows that $\cB$ is an integral $\cA$-algebra if and
only if $\cB_{!!}$ is an integral $\cA_{!!}$-algebra. Indeed,
suppose that $\cB$ is an integral $\cA$-algebra, so that
$\cB=\mathrm{i.c.}(\cA_{!!},\cB_{!!})^a$, whence by adjunction,
a morphism of $\cA_{!!}$-algebras
$\cB_{!!}\to\mathrm{i.c.}(\cA_{!!},\cB_{!!})$, which must be
the identity map.

(iii)\ \
Also, if $\cR''$ is the weak normalization of $\cR$ in $\cS$,
then it follows easily from proposition
\ref{prop_invariance-under-almiso}(i,ii) that $\cR''^a$ is the
weak normalization of $\cR^a$ in $\cS^a$ : the details are left
to the reader.
\end{remark}

\begin{proposition}\label{prop_integr-foll}
Suppose that the ideal $\fm$ fulfills condition $(\bB)$ of
\cite[\S2.1.6]{Ga-Ra}.
Let $X$ be a quasi-compact $S$-scheme, $\cA$ a quasi-coherent
$\cO_{\!X}^a$-algebra, and $\cB$ a flat $\cA$-algebra that is
almost finitely presented as an $\cA$-module. Then $\cB(X)$ is
an integral $\cA(X)$-algebra.
\end{proposition}
\begin{proof} For every $b\in\cB(X)_*$ and every $i\in\N$, set :
$$
c_i(b):=\tr_{\Lambda^i_\cA\cB/\cA}(\Lambda^i_\cA(b\one_\cB))
\in\cA(X)
$$
(see remark \ref{rem_global-trace-morph}(ii)). In view of
\cite[Lemma 8.2.28]{Ga-Ra}, it suffices to show

\begin{claim} For every $b\in\cB(X)_*$ and $a\in\fm$, there
exists $n\in\N$ such that :
\begin{enumerate}
\item
$c_i(ab)=0$ for every $i>n$.
\item
$\sum_{i=0}^n (-ab)^{i+1}\cdot c_{n-i}(ab)=0$.
\end{enumerate}
\end{claim}
\begin{pfclaim}[] Since $X$ admits a finite affine covering,
and since the trace morphism commutes with arbitrary base
changes, we are easily reduced to the case where $X$ is affine.
In this case, set $A:=\cA(X)$ and $B:=\cB(X)$; then $B$ is an
almost finitely generated projective $A$-algebra, and there
exists $n\in\N$ such that $a\one_B$ factors through $A$-linear
maps $u:B\to A^{\oplus n}$ and $v:A^{\oplus n}\to B$
(\cite[Lemma 2.4.15]{Ga-Ra}). Thus, $ab\one_B=v\circ(u\circ b\one_B)$,
and we get :
$$
c_i(ab)=
\tr_{\Lambda^i_AA^{\oplus n}/A}(\Lambda^i_A(u\circ b\one_B\circ v))
\qquad
\text{for every $i\in\N$}
$$
(\cite[Lemma 4.1.2(i)]{Ga-Ra}), from which (i) already follows.
By the same token, we also obtain :
$$
\chi:=\sum_{i=0}^n(-u\circ b\one_B\circ v)^i\cdot c_{n-i}(ab)=0
$$
(\cite[Prop.4.4.30]{Ga-Ra}). However, the left-hand side of the
identity of (ii) is none else than
$(-1)\cdot v\circ\chi\circ(u\circ b\one_B)$; whence the claim.
\end{pfclaim}
\end{proof}

The following lemma generalizes \cite[Cor. 4.4.31]{Ga-Ra}.

\begin{lemma}\label{lem_axiomatize-familiar}
Let $A$ be a $V^a$-algebra, $B\subset A$ a $V^a$-subalgebra,
$P$ an almost finitely generated projective $A$-module, and
$\phi$ an $A$-linear endomorphism of $P$. Suppose that :
\begin{enumerate}
\alphaenu
\item
$B$ is integrally closed in $A$.
\item
$\phi$ is integral over $B_*$.
\end{enumerate}
Then we have :
\begin{enumerate}
\item
$\tr_{P/A}(\phi)\in B_*$.
\item
If\/ $\fm$ fulfills condition $(\bB)$ of\/ \cite[\S2.1.6]{Ga-Ra},
then $\det(\one_P+T\phi)\in B_*[[T]]$.
\end{enumerate}
\end{lemma}
\begin{proof} For the meaning of assumption (b), see the proof
of \cite[Cor. 4.4.31]{Ga-Ra}. 

(i): For every $\eps\in\fm$ we may find a free $A$-module
$L$ of finite rank, and $A$-linear morphisms $u:P\to L$,
$v:L\to P$ such that $v\circ u=\eps\one_P$
(\cite[Lemma 2.4.15]{Ga-Ra}). It follows that 
$$
\tr_{L/A}(u\circ\phi\circ v)=\tr_{P/A}(\phi\circ v\circ u)=
\eps\cdot\tr_{P/A}(\phi)
$$
(\cite[Lemma 4.1.2]{Ga-Ra}). On the other hand, say that
the polynomial $T^n+\sum_{j=0}^{n-1}b_jT^j\in B_*[T]$ annihilates
$\phi$; it is easily seen that $T^n+\sum_{j=0}^{n-1}b_j\eps^{n-j}T^j$
annihilates $u\circ\phi\circ v$, so the latter is integral
over $B_*$ as well, and therefore its trace lies in $B_*$
(\cite[Cor.4.4.31 and Rem.8.2.30(i)]{Ga-Ra}).
Summing up, we have shown that $\eps\cdot\tr_{P/A}(\phi)\in B_*$
for every $\eps\in\fm$, whence the contention.

(ii): Recall that $\det(\one_P+T\phi)$ is the power series in
the variable $T$, whose coefficient in degree $i$ is the trace
of $\Lambda^i_A\phi$ on $\Lambda^i_AP$ : see
\cite[\S4.3.1, \S4.3.3]{Ga-Ra}.

\begin{claim}\label{cl_clever-trick}
Let $Q$ be another almost finitely generated projective
$A$-module, $\psi$ be an endomorphism of $P$ that is also
integral over $B_*$. Then the endomorphism $\phi\otimes_A\psi$
of $P\otimes_AQ$ is integral over $B_*$.
\end{claim}
\begin{pfclaim} The tensor product defines a map of unital
associative $B_*$-algebras
$$
\End_A(P)\otimes_{B_*}\End_A(Q)\to\End_A(P\otimes_AQ).
$$
By assumption, the subalgebra $B_*[\phi]\subset\End_A(P)$
is finite over $B_*$, and similarly for
$B_*[\psi]\subset\End_A(Q)$. Hence, the same holds for
the image of $B_*[\phi]\otimes_{B_*}B_*[\psi]$ in
$\End_A(P\otimes_AQ)$, whence the claim.
\end{pfclaim}

It follows easily from claim \ref{cl_clever-trick}, that
the endomorphism $\Lambda^i_A\phi$ of $\Lambda^i_AP$ is
integral over $B_*$. Hence, we may replace $P$ by
$\Lambda^i_AP$, and reduce to showing that the trace
of $\phi$ lies in $B_*$, which is known, by (i).
\end{proof}

\begin{lemma}\label{lem_detect-afr-after-quot}
Let $A$ be a $V^a$-algebra, $I\subset\rad(A)$ a tight ideal
(see \cite[Def.5.1.5]{Ga-Ra}), and $P$ a flat and almost
finitely generated $A$-module. Let also $r\in\N$. We have :
\begin{enumerate}
\item
If $P/IP$ is a faithfully flat $A$-module, then $P$ is a
faithfully flat $A$-module.
\item
If\/ $\fm$ fulfills condition $(\bB)$, and $P$ is an almost
projective $A$-module such that $P/IP$ is an $A/I$-module of
almost finite rank (resp. of finite rank $\leq r$), then $P$
is an $A$-module of almost finite rank (resp. of finite rank
$\leq r$).
\end{enumerate}
\end{lemma}
\begin{proof}(i): It suffices to show that for every almost
finitely generated $A$-module $M\neq 0$, we have
$M\otimes_AP\neq 0$. To this aim, it suffices to check that
$M\otimes_AP/IP\neq 0$; since
$M\otimes_AP/IP=(M\otimes_AA/I)\otimes_{A/I}P/IP$, and since
by assumption $P/IP$ is a faithfully flat $A/I$-module, we
are further reduced to showing that $M\otimes_AA/I\neq 0$;
the latter follows from \cite[Lemma 5.1.7]{Ga-Ra}.

(ii): Suppose first that $P/IP$ is of finite rank $\leq r$,
so that $\Lambda^r_{A/I}(P/IP)=0$. Hence
$(\Lambda^r_AP)\otimes_AP/IP=0$, and then
\cite[Lemma 5.1.7]{Ga-Ra} yields $\Lambda^r_AP=0$, as
sought. In the more general case where $P/IP$ has almost
finite rank, pick $n\in\N$ and a finitely generated
subideal $\fm_0\subset\fm$, such that $I^n\subset\fm_0$.
Let also $\fm_1\subset\fm$ be any finitely generated
subideal such that $\fm_0\subset\fm_1^{n+1}$; by assumption,
there exists $i\in\N$ such that $\fm_1\Lambda^i_{A/I}(P/IP)=0$,
{\em i.e.} $\fm_1\Lambda^i_AP\subset I\Lambda^i_AP$. Therefore
$$
\fm_0\Lambda^i_AP\subset\fm_1^{n+1}\Lambda^i_AP\subset
I^{n+1}\Lambda^i_AP\subset\fm_0I\Lambda^i_AP
$$
whence $\fm_0\Lambda^i_AP=0$, by \cite[Lemma 5.1.7]{Ga-Ra}.
Since $\fm_0$ is arbitrary, we deduce that $P$ has almost
finite rank, as claimed.
\end{proof}

\begin{lemma}\label{lem_detect-faith-after-inj-ext}
Let $A\to B$ be a morphism of $V^a$-algebras which is a
monomorphism on the underlying $V^a$-modules, $P$ an almost
finitely generated and almost projective $A$-module, and
$r\in\N$. Suppose as well that $\fm$ fulfills condition
$(\bB)$; the following holds :
\begin{enumerate}
\item
$P$ is an $A$-module of almost finite rank (resp. of finite
rank $\leq r$) if and only if the same holds for the
$B$-module $B\otimes_AP$.
\item
Suppose that $P$ has almost finite rank. Then $P$ is a
faithfully flat $A$-module if and only if the same holds
for the $B$-module $B\otimes_AP$.
\end{enumerate}
\end{lemma}
\begin{proof}(i): According to proposition
\ref{prop_ext-std-to-m-nonflat}(i), the $A$-module
$\Lambda^i_AP$ is flat for every $i\in\N$; hence the
natural $A$-linear morphism
$$
\Lambda^i_AP\to B\otimes_A\Lambda^i_AP\isom
\Lambda^i_B(B\otimes_AP)
$$
is a monomorphism for every $i\in\N$; the assertion follows
straightforwardly.

(ii): According to \cite[Th.4.3.28]{Ga-Ra} (which holds whenever
$\fm$ fulfills condition $(\bB)$) we have $V^a$-algebras
$A_0,A'$ and an isomorphism of $V^a$-algebras $A\isom A_0\times A'$
such that $A_0\otimes_AP=0$ and $A'\otimes_AP$ is a faithfully
flat $A'$-module. Suppose now that $B\otimes_AP$ is a faithfully
flat $B$-module; since
$(B\otimes_AA_0)\otimes_B(B\otimes_AP)\simeq B\otimes_AA_0\otimes_AP=0$,
it follows that $B\otimes_AA_0=0$. On the other hand, the induced
$A$-linear morphism $A_0\to B\otimes_AA_0$ is a monomorphism,
since $A_0$ is a flat $A$-module; thus, $A_0=0$, and consequently
$P=A'\otimes_AP$ is faithfully flat.
\end{proof}

\begin{lemma}\label{lem_extend-etale}
Let $j:U\to X$ be an open quasi-compact immersion of $S$-schemes,
$\cB$ a quasi-coherent $\cO^a_{\!U}$-algebra, $\cA$ a quasi-coherent
$\cO^a_{\!X}$-algebra, and $j^*\cA\to\cB$ an \'etale morphism.
Then the induced morphism
$$
\cA\to j_*j^*\cA\to j_*\cB
$$
is flat if and only if it is \'etale.
\end{lemma}
\begin{proof} We may assume that $j_*\cB$ is a flat
$\cA$-algebra, and it suffices to show that $j_*\cB$
admits a diagonal idempotent as in remark
\ref{rem_idemp-and-traces}(i). To this aim, let us
remark more generally :

\begin{claim}\label{cl_flat-obvious}
Let $f:Y\to X$ be a quasi-compact and quasi-separated morphism,
$\cF$ a flat quasi-coherent $\cA$-module, and $\cG$ an
$f^*\cA$-module, quasi-coherent as an $\cO^a_Y$-module.
Then the natural map
$$
\cF\otimes_\cA f_*\cG\to f_*(f^*\cF\otimes_{f^*\cA}\cG)
$$
is an isomorphism.
\end{claim}
\begin{pfclaim} This is proved just as for $\cO_{\!X}$-modules.
Namely, one reduces easily to the case where $X$ is affine,
and one needs to check that the natural map
$$
\cF(X)\otimes_{\cA(X)}\cG(Y)\to
(f^*\cF\otimes_{f^*\cA}\cG)(Y)
$$
is an isomorphism. However, by assumption $Y$ can be covered
by finitely many affine open subsets $U_1,\dots,U_n$, and the
intersection $U_{ij}:=U_i\cap U_j$ can be covered by finitely
many affine open subsets $U_{ij1},\dots,U_{ijn}$, for every
$i,j=1,\dots,n$. Since $f^*\cF$, $\cG$ and $f^*\cA$
are quasi-coherent $\cO^a_Y$-modules, the natural maps
$$
\begin{aligned}
\cF(X)\otimes_{\cA(X)}\cG(U_i)\,&\to
(f^*\cF\otimes_{f^*\cA}\cG)(U_i) \\
\cF(X)\otimes_{\cA(X)}\cG(U_{ijk})\,&\to
(f^*\cF\otimes_{f^*\cA}\cG)(U_{ijk})
\end{aligned}
$$
are isomorphisms, for every $i,j,k=1,\dots,n$. Then we get :
$$
\begin{aligned}
(f^*\cF\otimes_\cA\cG)(Y)=\, &
\xymatrix{\Ker(\prod_{i=1}^n\cF(X)\otimes_{\cA(X)}\cG(U_i)
\ar@<.5ex>[r] \ar@<-.5ex>[r] &
\prod_{i,j,k=1}^n\cF(X)\otimes_{\cA(X)}\cG(U_{ijk}))} \\
=\, &
\cF(X)\otimes_{\cA(X)}\xymatrix{\Ker(\prod_{i=1}^n\cG(U_i)
\ar@<.5ex>[r] \ar@<-.5ex>[r] & \prod_{i,j,k=1}^n\cG(U_{ijk}))} 
\end{aligned}
$$
(since $\cF(X)$ is a flat $\cA(X)$-module) whence the claim.
\end{pfclaim}

Now, set $\cR:=\cB\otimes_{j^*\cA}\cB$; since $j$ is
quasi-compact and $j_*\cB$ is a flat $\cA$-algebra,
claim \ref{cl_flat-obvious} implies that the natural
morphism :
$$
j_*\cB\otimes_{\cA}j_*\cB\to j_*\cR
$$
is an isomorphism. Especially, the diagonal idempotent
of $\cB$ extends to a global section of
$j_*\cB\otimes_{\cA}j_*\cB$, and the assertion follows.
\end{proof}

The criterion of lemma \ref{lem_extend-etale} has
limited usefulness, since it is not usually known
a priori that $j_*\cB$ is a flat $\cA$-algebra.
In several situations, one can however apply the
following variant.

\begin{proposition}\label{prop_axiomatize-familiar}
Let $X$ be an $S$-scheme, $j:U\to X$ a quasi-compact open
immersion, $\cA\to\cB$ a morphism of quasi-coherent
$\cO^a_{\!X}$-algebras, and suppose that :
\begin{enumerate}
\alphaenu
\item
The units of adjunction $\cA\to j_*j^*\cA$ and $\cB\to j_*j^*\cB$
are monomorphisms.
\item
$\cA$ is integrally closed in $j_*j^*\cA$, and $\cB$ is an
integral $\cA$-algebra.
\item
$j^*\cB$ is an \'etale $j^*\cA$-algebra and an almost finitely
presented $j^*\cA$-module.
\item
The diagonal idempotent $e_{j^*\cB/j^*\cA}$ is a global section
of the subalgebra
$$
\Img(\cB\otimes_\cA\cB\to j_*j^*(\cB\otimes_\cA\cB)).
$$
\end{enumerate}
Then $\cB$ is an \'etale $\cA$-algebra and an almost finitely
presented $\cA$-module.
\end{proposition}
\begin{proof} We may assume that $X$ is affine.
Under our assumptions, the restriction map $\cA(X)_*\to\cA(U)_*$
is injective, and its image is integrally closed in $\cA(U)_*$
(\cite[Rem.8.2.30]{Ga-Ra}); moreover, $\cB(X)_{!!}$ is an integral
$\cA(X)_{!!}$-algebra (remark \ref{rem_alm-int-closure}(ii)).
In view of lemma \ref{lem_axiomatize-familiar}(i), we deduce a
commutative diagram
\set\begin{equation}\label{eq_nausee}
{\diagram
\cB\otimes_\cA\cB \ar[rrr]^-t \ar[d] & & & \cA \ar[d] \\
j_*j^*(\cB\otimes_\cA\cB) \ar[rrr]^-{j_*(t_{j^*\cB/j^*\cA})}
& & & j_*j^*\cA
\enddiagram}
\end{equation}
whose vertical arrows are the units of adjunctions, and
where $t$ is a uniquely determined $\cA$-bilinear form.
Then we are reduced to showing :

\begin{claim}\label{cl_nausee}
Let $j:U\to X$ be as in the proposition, and $\cA\to\cB$ a
morphism of quasi-coherent $\cO^a_{\!X}$-algebras fulfilling
conditions (a), (c), (d) of the proposition, and such that
there exists a bilinear form $t$ making commute diagram
\ref{eq_nausee}. Then $\cB$ is an \'etale $\cA$-algebra
and an almost finitely presented $\cA$-module.
\end{claim}
\begin{pfclaim}[] The proof proceeds by a familiar argument
(see {\em e.g.} the proof of \cite[Claim 3.5.33]{Ga-Ra}).
Namely, for a given $\eps\in\fm$, we can write
$$
\eps\cdot e_{j^*\cB/j^*\cA}=\sum_{j=1}^m a_j\otimes b_j
\qquad
\text{where $a_j,b_j\in\cB(X)_*$ for every $j=1,\dots,m$.}
$$
In light of remark \ref{rem_idemp-and-traces}(iv), we deduce
$$
\eps\cdot a=\sum_{j=1}^m t(a\otimes a_j)\cdot b_j
\qquad
\text{for every $a\in\cB(X)_*$.}
$$
Indeed, the identity holds after restriction to $U$,  and
assumption (a) implies that the restriction map
$\cB(X)_*\to\cB(U)_*$ is injective. We may then define
$\cA$-linear maps :
$\cB\xrightarrow{\phi}\cA^{\oplus m}\xrightarrow{\psi}\cB$
by the rules :
$$
\phi(a):=(t(a\otimes a_1),\dots,t(a\otimes a_m)) \qquad
\psi(s_1,\dots,s_m):=\sum^m_{j=1}s_jb_j
$$
for every open subset $U\subset X$, every $a\in\cB(U)_*$,
and every $s_1,\dots,s_m\in\cA(U)_*$. Thus we have
$\psi\circ\phi=\eps\one_\cB$, and since $\eps$ is arbitrary,
this already proves that $\cB$ is an $\cA$-module of almost
finite type, so we may apply
\cite[Lemma 2.4.15 and Prop.2.4.18(i)]{Ga-Ra} to deduce that
$\cB$ is a flat and almost finitely presented $\cA$-module.
Now, assumption (a) and claim \ref{cl_flat-obvious} imply
that the horizontal arrows of the induced diagram
$$
\xymatrix{ \cB\otimes_\cA\cB \ar[r] \ar[d]_{\mu_{\cB/\cA}} &
j_*j^*(\cB\otimes_\cA\cB) \ar[d]^{j_*j^*\mu_{\cB/\cA}} \\
\cB \ar[r] & j_*j^*\cB
}$$
are monomorphisms. From the characterization of remark
\ref{rem_idemp-and-traces}(i), it follows easily that
the section $e_{j^*\cB/j^*\cA}$ is the diagonal idempotent
for the morphism $\cA\to\cB$, so $\cB$ is an unramified
$\cA$-algebra, and the proof is concluded.
\end{pfclaim}
\end{proof}

\begin{corollary}\label{cor_axiomatize-familiar}
Let $X$ be an $S$-scheme, $j:U\to X$ a quasi-compact open
immersion, $\cA\to\cB$ a morphism of quasi-coherent
$\cO^a_{\!U}$-algebras, and suppose that :
\begin{enumerate}
\alphaenu
\item
$\cB$ is an \'etale $\cA$-algebra and an almost finitely
presented $\cA$-module.
\item
The diagonal idempotent $e_{\cB/\cA}$ is a global section of
the subalgebra
$$
\Img((j_*\cB)\otimes_{j_*\cA}(j_*\cB)\to j_*(\cB\otimes_\cA\cB)).
$$
\end{enumerate}
Then $j_*\cB$ is an \'etale $j_*\cA$-algebra and an almost
finitely presented $j_*\cA$-module.
\end{corollary}
\begin{proof} Indeed, the induced morphism $j_*\cA\to j_*\cB$
fulfills conditions (a), (c) and (d) of proposition
\ref{prop_axiomatize-familiar}, and diagram \ref{eq_nausee}
trivially commutes with $t:=j_*(t_{\cB/\cA})$, so the
corollary follow from claim \ref{cl_nausee}.
\end{proof}

Using lemma \ref{lem_axiomatize-familiar} we can also
relax one assumption in \cite[Prop.8.2.31(i)]{Ga-Ra};
namely, we have the following :

\begin{proposition}\label{prop_8231}
Let $A\subset B$ be a pair of $V^a$-algebras, such that
$A=\mathrm{i.c.}(A,B)$. Then, for every \'etale almost
finite projective $A$-algebra $A_1$ we have
$A_1=\mathrm{i.c.}(A_1,A_1\otimes_AB)$.
\end{proposition}
\begin{proof} Set $B_1:=A_1\otimes_AB$, and suppose that
$x\in B_{1*}$ is integral over $A_{1*}$. Let $e\in(A_1\otimes_AA_1)_*$
be the diagonal idempotent of the unramified $A$-algebra $A_1$;
for given $\eps_1,\eps_2,\eps_3\in\fm$ we write
$\eps_1\cdot e=\sum^k_{j=1}c_i\otimes d_i$ for some 
$c_i,d_i\in A_{1*}$. According to remark
\ref{rem_idemp-and-traces}(iv) and \cite[Prop.4.1.8(ii)]{Ga-Ra},
we have $\sum^k_{j=1}c_i\cdot\Tr_{B_1/B}(xd_i)=\eps\cdot x$.
On the other hand, $\eps_2\cdot x$ and $\eps_3\cdot d_i$ are
integral over $A_*$ (by \cite[Lemma 5.1.13(i)]{Ga-Ra}, which
holds for arbitrary basic setups); then lemma
\ref{lem_axiomatize-familiar}(i) implies that
$\Tr_{B_1/B}(\eps_2\eps_3\cdot xd_i)\in A_*$ for every $i=1,\dots,k$.
Hence $\eps_1\eps_2\eps_3\cdot x\in A_{1*}$, and since
$\eps_1,\eps_2,\eps_3$ are arbitrary, the assertion follows. 
\end{proof}

\begin{lemma}\label{lem_complem-qcoh-alg}
Let $X$ be an $S$-scheme, $\cA$ a quasi-coherent
$\cO^a_{\!X}$-algebra. We have :
\begin{enumerate}
\item
Let $f:Y\to X$ a faithfully flat and quasi-compact morphism.
Then $\cA$ is an \'etale (resp. weakly \'etale, resp. weakly
unramified) $\cO^a_{\!X}$-algebra if and only if $f^*\cA$
is an \'etale (resp. weakly \'etale, resp. weakly unramified)
$\cO^a_Y$-algebra.
\item
Suppose that $X$ is integral and $\cA$ is a weakly \'etale
$\cO^a_{\!X}$-algebra. Let $\eta\in X$ be the generic point;
then $\cA$ is an \'etale $\cO^a_{\!X}$-algebra if and only if
$\cA_\eta$ is an \'etale $\cO^a_{\!X,\eta}$-algebra.
\item
Suppose that $X$ is normal and irreducible, and $\cA$ is
integral, torsion-free and unramified over $\cO^a_{\!X}$.
Then $\cA$ is {\'e}tale over $\cO^a_{\!X}$, and almost
finitely presented as an $\cO^a_{\!X}$-module.
\item
Let $j:X\to Y$ be an open quasi-compact immersion of
$S$-schemes, with $Y$ normal and irreducible, and
such that $j_*\cO_{\!X}=\cO_{\!Y}$, and suppose that
$\cA$ is integral and torsion-free over $\cO^a_{\!X}$.
Then $j_*\cA$ is an integral $\cO^a_{\!Y}$-algebra.
\end{enumerate}
\end{lemma}
\begin{proof}(i): Arguing as in the proof of lemma
\ref{lem_complem-qcoh}, one reduces easily to the case where
$X$ and $Y$ are affine, in which case the assertion follows
from \cite[\S3.4.1]{Ga-Ra}.

(ii): Let $U\subset X$ be any affine open subset; by assumption
$\cA(U)$ is flat over the almost ring
$B:=\cA(U)\otimes_{\cO^a_{\!X}(U)}\cA(U)$, and
$\cA(U)\otimes_{\cO^a_{\!X}(U)}\kappa(\eta)^a$ is almost projective
over $B\otimes_{\cO^a_{\!X}(U)}\kappa(\eta)^a$. Then the assertion
follows from \cite[Prop.2.4.19]{Ga-Ra}.

(iii): Let $\eta\in X$ be the generic point. We begin with the
following :

\begin{claim}\label{cl_generically-true}
Suppose that $V=\fm$ (the ``classical limit'' case of
\cite[Ex.2.1.2(ii)]{Ga-Ra}), and let $f:A\to B$ be
a local homomorphism of local rings. Then :
\begin{enumerate}
\item
$f^a$ is weakly {\'e}tale if and only if $f$ extends to
an isomorphism of strict henselizations $A^\sh\to B^\sh$.
\item
Especially, if $A$ is a field and $f^a$ is weakly {\'e}tale,
then $B$ is a separable algebraic extension of $A$.
\end{enumerate}
\end{claim}
\begin{pfclaim} (i): In the classical limit case, a weakly
{\'e}tale morphism is the same as an ``absolutely flat'' map
as defined in \cite{Ol}; then the assertion is none else
than \cite[Th.5.2]{Ol}.

(ii): If $A$ is a field, $A^\sh$ is a separable closure of
$A$; hence (i) implies that $B$ is a subring of $A^\sh$, hence
it is a subfield of $A^\sh$.
\end{pfclaim}

\begin{claim}\label{cl_generically-truebis}
Suppose that $V=\fm$. Let $A$ be a field and $f:A\to B$ a ring
homomorphism such that $f^a$ is weakly {\'e}tale. Then :
\begin{itemize}
\item[]
\begin{enumerate}
\item
Every finitely generated $A$-subalgebra of $B$ is finite {\'e}tale
over $A$ (that is, in the usual sense of \cite[Ch.IV, \S17]{EGA4}).
\item
$f^a$ is {\'e}tale if and only if $f$ is {\'e}tale (in the
usual sense of \cite{EGA4}).
\end{enumerate}
\end{itemize}
\end{claim}
\begin{pfclaim} (i): From claim \ref{cl_generically-true}(ii) we see
that $B$ is reduced of Krull dimension $\leq 0$, and all its residue
fields are separable algebraic extensions of $A$. We consider first
the case of a monogenic extension $C:=A[b]\subset B$.
For every prime ideal $\fp\subset B$, let $b_\fp$ be the image of
$b$ in $B_\fp$; then for every such $\fp$ we may find an irreducible
separable monic polynomial $P_\fp(T)\in A[T]$ with $P_\fp(b_\fp)=0$.
This identity persists in an open neighborhood $U_\fp\subset\Spec\,B$
of $\fp$, and finitely many such $U_\fp$ suffice to cover $\Spec\,B$.
Thus, we find finitely many $P_{\fp_1}(T),\dots,P_{\fp_k}(T)$
such that $\prod_{i=1}^kP_{\fp_i}(b)=0$ holds in $B$, and after
omitting repetitions, we may assume that all these polynomials
are distinct. Since $C$ is a quotient of the separable $A$-algebra
$A[T]/(\prod_{i=1}^kP_{\fp_i}(T))$, the claim follows in this case.
In the general case, we may write $C=A[b_1,\dots,b_n]$ for certain
$b_1,\dots,b_n\in B$. Then $\Spec\,C$ is a reduced closed subscheme of
$\Spec\,A[b_1]\times_{\Spec\,A}\cdots\times_{\Spec\,A}\Spec\,A[b_n]$;
by the foregoing, the latter is \'etale over $\Spec\,A$, hence the
same holds for $\Spec\,C$.

(ii): We may assume that $f^a$ is \'etale, and we have to show
that $f$ is \'etale, {\em i.e.} that $B$ is a finitely generated
$A$-module. Hence, let $e_{B/A}\in B\otimes_AB$ be the diagonal
idempotent (see remark \ref{rem_idemp-and-traces}(i)); we choose
a finitely generated $A$-subalgebra $C\subset B$ such that $e_{B/A}$
is the image of an element $e'\in C\otimes_AC$. Notice that $1-e'$
lies in the kernel of the multiplication map
$\mu_{C/A}:C\otimes_AC\to C$. By (i), the morphism $A\to C$
is \'etale, hence it admits as well a diagonal idempotent
$e_{C/A}\in C\otimes_AC$; on the other hand,
\cite[Lemma 3.1.2(v)]{Ga-Ra} says that $B^a$ is \'etale over
$C^a$, especially $B$ is a flat $C$-algebra, hence the natural
map $C\otimes_AC\to B\otimes_AB$ is injective.
Since $1-e_{C/A}$ lies in the kernel of the multiplication map
$\mu_{B/A}:B\otimes_AB\to B$, we have :
$$
e_{B/A}(1-e_{C/A})=0=e_{C/A}(1-e')=e_{C/A}(1-e_{B/A})
$$
from which it follows easily that $e_{B/A}=e_{C/A}$.
Moreover, the induced morphism $\Spec\,B\to\Spec\,C$ has dense
image, so it must be surjective, since $\Spec\,C$ is a discrete
finite set; therefore $B$ is even a faithfully flat $C$-algebra.
Let $J:=\Ker(B\otimes_AB\to B\otimes_CB)$; then $J$ is the ideal
generated by all elements of the form $1\otimes c-c\otimes 1$,
where $c\in C$; clearly this is the same as the extension of the
ideal $I_{C/A}:=\Ker\,\mu_{C/A}$. However, $I_{C/A}$
is generated by the idempotent $1-e_{C/A}$ (\cite[Cor.3.1.9]{Ga-Ra}),
consequently $J$ is generated by $1-e_{B/A}$, {\em i.e.}
$J=\Ker\,\mu_{B/A}$.
So finally, the multiplication map $\mu_{B/C}:B\otimes_CB\to B$ is an
isomorphism, whose inverse is the map $B\to B\otimes_CB$ :
$b\mapsto b\otimes 1$. The latter is of the form $j\otimes_C\one_B$,
where $j:C\to B$ is the natural inclusion map. By faithfully flat
descent we conclude that $C=B$, whence the claim.
\end{pfclaim}

Next, we remark that $\cA_\eta$ is a finitely presented
$\cO^a_{\!X,\eta}$-module. Indeed, since $K:=\cO_{\!X,\eta}$
is a field, we have either $\fm K=0$, in which case the category
of $K^a$-modules is trivial and there is nothing to show, or
else $\fm K=K$. So we may assume that we are in the
``classical limit'' case, and then the assertion
follows from claim \ref{cl_generically-truebis}(ii).

Now, set $\cB:=\cA\otimes_{\cO^a_{\!X}}\cA$, and
let $j:X(\eta)\to X$ be the natural morphism; since
$X$ is normal and $\cA$ is torsion-free, the units
of adjunction $\cO^a_{\!X}\to j_*j^*\cO_{\!X(\eta)}$
and $\cA\to j_*j^*\cA$ are monomorphisms. Since $\cA$
is unramified over $\cO^a_{\!X}$, the diagonal idempotent
of $j^*\cA$ lies in the image of the restriction map
$\cB(X)_*\to\cB_{\eta*}$. In view of these observations,
an easy inspection shows that the proof of proposition
\ref{prop_axiomatize-familiar} carries over {\em verbatim\/}
to the current situation, and yields assertion (iii).

(iv): Let $\cT\subset\cA_{!!}$ be the maximal
torsion $\cO_{\!X}$-subsheaf, and set $\cB:=\cA_{!!}/\cT$.
Then $\cB$ is an integral, torsion-free $\cO_{\!X}$-algebra,
by remark \ref{rem_alm-int-closure}(ii), and $\cB^a\simeq\cA$,
hence $(j_*\cB)^a\simeq j_*\cA$. Let $\eta\in Y$ be the generic
point; in light of remark \ref{rem_alm-int-closure}(ii), it
then suffices to show :

\begin{claim} Under the assumptions of (v), let $\cR$ be
an integral quasi-coherent and torsion-free $\cO_{\!X}$-algebra.
Then $j_*\cR$ is an integral $\cO_{\!Y}$-algebra.
\end{claim}
\begin{pfclaim}[] Let us write $\cR_\eta$ as the filtered
union of the family $(R_\lambda~|~\lambda\in\Lambda)$ of
its finite $\cO_{\!Y,\eta}$-subalgebras, and for every
$\lambda\in\Lambda$, let $\cR_\lambda\subset\cR$ be the
quasi-coherent $\cO_{\!X}$-subalgebra such that
$\cR_\lambda(V)=R_\lambda\cap\cR(V)$ for every non-empty open
subset $V\subset X$. Then $\cR$ is the filtered colimit of the
family $(\cR_\lambda~|~\lambda\in\Lambda)$,
and clearly it suffices to show the claim for every $\cR_\lambda$.
We may thus assume from start that $\cR_\eta$ is a finite
$\cO_{Y,\eta}$-algebra. Let $\cR^\nu$ be the integral closure
of $\cR$, {\em i.e.} the quasi-coherent $\cO_{\!X}$-algebra
such that $\cR^\nu(V)$ is the integral closure of $\cR(V)$
in $\cR_\eta$, for every non empty affine open subset $V\subset X$.
Clearly it suffices to show that $j_*\cR^\nu$ is integral,
hence we may replace $\cR$ by $\cR^\nu$ and assume from
start that $\cR$ is integrally closed. In this case, for
every non-empty open subset $V\subset X$ the restriction map
$\cR(V)\to\cR_\eta$ induces a bijection between
the idempotents of $\cR(V)$ and those of $\cR_\eta$.
Especially, $\cR(X)$ admits finitely many idempotents,
and moreover we have a natural decomposition
$\cR=\cR_1\times\cdots\times\cR_k$, as a product of
$\cO_{\!X}$-algebras, such that $\cR_{i,\eta}$ is a
field for every $i\leq k$. It then suffices to show
the claim for every $j_*\cR_i$, and then we may assume
throughout that $\cR_\eta$ is a field. Up to replacing
$\cR$ by its normalization in a finite extension of
$\cR_\eta$, we may even assume that $\cR_\eta$ is
a finite normal extension of $\cO_{\!Y,\eta}$.
Hence, let $V\subset Y$ be any non-empty affine open
subset, and $a\in\cR(X\cap V)$ any element. Let $P(T)$
be the minimal polynomial of $a$ over the field $\cO_{Y,\eta}$;
we have to show that the coefficients of $P(T)$ lie
in $\cO_{\!Y}(V)=\cO_{\!X}(X\cap V)$, and since $Y$
is normal, it suffices to prove that these coefficients
are integral over $\cO_{\!X}(W)$, for every non-empty
affine open subset $W\subset X\cap V$. However, since
$\cR_\eta$ is normal over $\cO_{Y,\eta}$, such coefficients
can be written as some elementary symmetric polynomials of the
conjugates of $a$ in $\cR_\eta$. Hence, we come down to showing
that the conjugates of $a$ are still integral over
$\cO_{\!X}(W)$. The latter assertion is clear: indeed,
if $Q(T)\in\cO_{\!X}(W)[T]$ is a monic polynomial
with $Q(a)=0$, then we have also $Q(a')=0$ for every
conjugate $a'$ of $a$.
\end{pfclaim}
\end{proof}

\begin{proposition}\label{prop_approx-al-alg}
Let $X$ be a quasi-compact and quasi-separated $S$-scheme,
$j:U\to X$ a quasi-compact open immersion, and $\cA$ a
quasi-coherent $\cO^a_{\!U}$-algebra, almost finitely
presented as an $\cO^a_{\!U}$-module.
Then for every finitely generated subideal $\fm_0\subset\fm$, there
exist a quasi-coherent $\cO_{\!X}$-algebra $\cB$, finitely presented
as an $\cO_{\!X}$-module, and a morphism $\cB^a_{|U}\to\cA$ of
$\cO^a_{\!U}$-algebras, whose kernel and cokernel are annihilated
by $\fm_0$.
\end{proposition}
\begin{proof} Set $\cC:=\mathrm{i.c.}(\cO^a_{\!X},j_*\cA)$;
then $\cC$ is a quasi-coherent $\cO^a_{\!X}$-algebra, and
$\cC_{|U}=\cA$, by lemma \ref{lem_about-integral-closure}.
According to proposition \ref{prop_fp-approx}, we may write
$\cC_{!!}$ as the colimit of a filtered family $(\cF_i~|~i\in I)$
of finitely presented quasi-coherent $\cO_{\!X}$-modules.
Pick a finitely generated subideal $\fm_1\subset\fm$ such that
$\fm_0\subset\fm_1^2$; for every affine open subset $U'\subset U$,
we may find $i\in I$ such that
$\fm_1\cdot\cA_{!!}(U')\subset\Img(\cF_i(U')\to\cA_{!!}(U'))$,
and since $U$ is quasi-compact, finitely many such open subsets
cover $U$; hence, we may find $i\in I$ such that
$\fm_1\cdot\cA_{!!}\subset\Img(\cF_{i|U}\to\cA_{!!})$; we set
$\cF:=\cF_i$, and let $\phi:\cG:=\Sym^\bullet_{\cO_{\!X}}\cF\to\cC_{!!}$
be the induced morphism of quasi-coherent $\cO_{\!X}$-algebras.
Notice that $\cG$ is finitely presented as an $\cO_{\!X}$-algebra.
Using again proposition \ref{prop_fp-approx} (or
\cite[Ch.I, Cor.9.4.9]{EGAI}) we may write $\Ker\,\phi$ as the
colimit of a filtered family $(\cK_\lambda~|~\lambda\in\Lambda)$
of quasi-coherent $\cO_{\!X}$-submodules of finite type.
Fix a finite covering $(U'_i~|~i\in I)$ of $X$ consisting
of affine open subsets, and for every $i\in I$, let
$f_{i,1},\dots,f_{i,n(i)}$ be a finite set of generators
of $\cF(U'_i)$; we may find monic polynomials
$P_{i,1}(T),\dots,P_{i,n(i)}(T)$ with coefficients in
$\cO_{\!X}(U'_i)$ such that $P_{i,j}(\phi(f_{i,j}))=0$ in
$\cC(U'_i)_{!!}$, for every $i\in I$ and every $j\leq n(i)$.
Let $\lambda\in\Lambda$ be chosen large enough, so that
$P_{i,j}(f_{i,j})\in\cK_\lambda(U'_i)$ for every $i\in I$
and every $j\leq n(i)$. Let also $\cK'_\lambda\subset\cG$
be the ideal generated by $\cK_\lambda$; then
$\cG':=\cG/\cK'_\lambda$ is an integral $\cO_{\!X}$-algebra of
finite presentation, hence it is finitely presented as an
$\cO_{\!X}$-module, in view of claim \ref{cl_fin-present}.
Clearly $\phi_{|U}$ descends to a map $\phi':\cG'_{|U}\to\cA_{!!}$,
whose cokernel is annihilated by $\fm_1$, and by
\cite[Claim 2.3.12 and the proof of Cor.2.3.13]{Ga-Ra} we may
find, for every affine open subset $U'\subset U$, a finitely
generated submodule
$K_{U'}\subset\Ker\,(\phi'_{U'}:\cG'(U')\to\cA(U')_{!!})$
such that $\fm_1^2\cdot\Ker\,\phi'_{U'}\subset K_{U'}$.
Another invocation of proposition \ref{prop_fp-approx}
ensures the existence of a quasi-coherent $\cO_{\!X}$-submodule
of finite type $\cJ\subset\Ker\,\phi$ such that $K_{U'}\subset\cJ(U')$
for all the $U'$ of a finite covering of $U$; the $\cO_{\!X}$-algebra
$\cB:=\cG'/\cJ\cdot\cG'$ fulfills the required conditions.
\end{proof}

\sset\subsubsection{}\label{subsec_approx-etale}
Let $(X_i~|~i\in I)$ be a cofiltered system of quasi-compact and
quasi-separated $S$-schemes, with affine transition morphisms
$h_\phi:X_j\to X_i$, for every morphism $\phi:j\to i$ in $I$.
Let also $U_i\subset X_i$ be a quasi-compact open subset, for
every $i\in I$, such that $U_j=h_\phi^{-1}U_i$ for every
$\phi:j\to i$ in $I$. Set
$$
X:=\lim_{i\in I}X_i
\qquad
U:=\lim_{i\in I}U_i.
$$
and denote by $h_i:U\to X_i$ the natural morphism, for
every $i\in I$.

\begin{corollary}\label{cor_approx-etale}
In the situation of \eqref{subsec_approx-etale}, let $\cA$ be any
flat almost finitely presented $\cO^a_{\!U}$-algebra. Then, for
every finitely generated subideal $\fm_0\subset\fm$ there exist
$i\in I$, a quasi-coherent $\cO_{\!X_i}$-algebra $\cR$, finitely
presented as an $\cO_{\!X_i}$-module, and a map of
$\cO^a_{\!U}$-algebras $f:h^*_i\cR^a\to\cA$ such that :
\begin{enumerate}
\alphaenu
\item
$\Ker\,f$ and $\Coker\,f$ are annihilated by $\fm_0$.
\item
For every $x\in U_i$ and every $b\in\fm_0$, the map
$b\cdot\one_{\cR,x}:\cR_x\to\cR_x$ factors through a
free $\cO_{\!U_i,x}$-module.
\end{enumerate}
\end{corollary}
\begin{proof} From \cite[Ch.IV, Th.8.3.11]{EGAIV-3} it is easily
seen that $X$ is a quasi-compact and quasi-separated $S$-scheme,
and $U$ is a quasi-compact open subset of $X$. Pick a finitely
generated subideal $\fm_1$ such that $\fm_0\in\fm_1^2$; according to
proposition \ref{prop_approx-al-alg}, we may find a quasi-coherent
$\cO_{\!X}$-algebra $\cB$, finitely presented as an $\cO_{\!X}$-module,
and a morphism $g:\cB^a_{|U}\to\cA$ of $\cO_{\!U}^a$-algebras whose
kernel and cokernel are annihilated by $\fm_1$. By standard
arguments we deduce
\set\begin{equation}\label{eq_tor-of-B}
\fm^2_1\cdot\Tor^1_{\cO_X(U')}(M,\cB(U'))=0
\end{equation}
for every affine open subset $U'\subset U$ and every
$\cO_{\!X}(U')$-module $M$. Now, let us remark, quite
generally :

\begin{claim}\label{cl_never-before}
Let $A$ be any ring, $M$ and $N$ any two $A$-modules,
$\phi:M\to N$ any $A$-linear map, and $M^\vee:=\Hom_A(M,A)$.
Then $\phi$ factors through a free $A$-module of finite rank
if and only if it lies in the image of the natural map
\set\begin{equation}\label{eq_never-before}
M^\vee\otimes_AN\to\Hom_A(M,N).
\end{equation}
\end{claim}
\begin{pfclaim} Suppose first that $\phi$ factors as the
composition of $A$-linear maps $\psi:M\to A^{\oplus n}$
and $\psi':A^{\oplus n}\to N$, for some $n\in\N$. For every
$i=1,\dots,n$, let $p_i:A^{\oplus n}\to A$ (resp.
$e_i:A\to A^{\oplus n}$) be the natural epimorphism
(resp. monomorphism), and set $\psi_i:=p_i\circ\psi$,
$\psi'_i:=\psi'\circ e_i$. Then
$\phi=\sum_{i=1}^n\psi'_i\circ\psi_i$ and each
summand $\psi'_i\circ\psi_i$ lies in the image
of \eqref{eq_never-before}. Conversely, if
$\phi$ lies in the image of \eqref{eq_never-before},
we have $\phi=\sum_{i=1}^n\psi'_i\circ\psi_i$, for
some $\psi_1,\dots,\psi_n\in M^\vee$ and
$\psi'_1,\dots,\psi'_n\in\Hom_A(A,N)=N$, and
we may define $\psi:=\sum_{i=1}^ne_i\circ\psi_i$,
$\psi'_i:=\sum_{i=1}^n\psi'_i\circ p_i$. Then
$\phi=\psi'\circ\psi$.
\end{pfclaim}

Now, pick a finite covering
$(U'_{i,\lambda}~|~\lambda\in\Lambda)$ of $U_i$,
consisting of affine open subsets, and a finite
system $b_1,\dots,b_k$ of generators of $\fm_1^2$,
and set $U'_\lambda:=h_i^{-1}U'_{i,\lambda}$ for every
$\lambda\in\Lambda$; from claim \ref{cl_never-before},
\eqref{eq_tor-of-B} and \cite[Lemma 2.4.17]{Ga-Ra}
we deduce that for every $\lambda\in\Lambda$ and
every $j=1,\dots,k$ there exists $n(\lambda,j)\in\N$
such that $b_j\cdot(\one_\cB)_{|U'_\lambda}$ factors through
$\cO^{\oplus n(\lambda,j)}_{\!U'_\lambda}$. Next, by
virtue of \cite[Ch.IV, Th.8.5.2]{EGAIV-3} we may
assume that $\cB$ descends to a $\cO_{\!X_i}$-algebra
$\cR$, finitely presented as an $\cO_{\!X_i}$-module, and
such that $b_j\cdot(\one_{\cR})_{|U'_{i,\lambda}}$ factors
through $\cO_{\!U'_{i,\lambda}}^{\oplus n(\lambda,j)}$ for
every $j=1,\dots,k$ and every $\lambda\in\Lambda$.
It follows easily that
$$
\fm_1^2\cdot\Tor^1_{\cO_{X_i}(U'_{i,\lambda})}(M,\cR(U'_{i,\lambda}))=0
$$
for every $\lambda\in\Lambda$ and every
$\cO_{X_i}(U'_{i,\lambda})$-module $M$. Then,
again \cite[Lemma 2.4.17]{Ga-Ra} and claim
\ref{cl_never-before} imply that condition
(b) holds for $\cR$, and by construction we also
get condition (a).
\end{proof}

\subsection{Almost pure pairs}\label{subsec_almost-pure-pair}
Throughout this section, we keep the notation of
\eqref{subsec_qcoh-almost-rings}, and {\em we assume
that $\fm$ fulfills condition $(\bB)$}.

\begin{definition}\label{def_al-pure}
Let $X$ be an $S$-scheme, $Z\subset X$ a closed subscheme
such that $U:=X\!\setminus\!Z$ is a dense subset of $X$,
and denote by $j:U\to X$ the open immersion.
\begin{enumerate}
\item
We say that the pair $(X,Z)$ is {\em almost pure relative to
$(V,\fm)$\/} (or just {\em almost pure}, when the underlying
basic setup is clear from the context) if the restriction
functor
\set\begin{equation}\label{eq_resr_etale}
\cO^a_{\!X}\Et_\mathrm{afr}\to\cO^a_{\!U}\Et_\mathrm{afr}
\quad :\quad \cA\mapsto\cA_{|U}
\end{equation}
from the category of \'etale $\cO^a_{\!X}$-algebras of almost
finite rank, to the category of \'etale $\cO^a_{\!U}$-algebras
of almost finite rank, is an equivalence.
\item
We say that the pair $(X,Z)$ is {\em normal\/} if $Z$
is a constructible subset of $X$, and the natural map
$\cO^a_{\!X}\to j_*\cO^a_{\!U}$ is a monomorphism, whose
image is integrally closed in $j_*\cO^a_{\!U}$.
\end{enumerate}
\end{definition}

\begin{remark}\label{rem_almost-pairs}
Let $(X,Z)$ be a pair as in definition \ref{def_al-pure},
where $Z$ is a constructible subset of $X$, and set
$U:=X\!\setminus\!Z$.

(i)\ \
Consider the following conditions:
\begin{enumerate}
\alphaenu
\item
$\cO_{\!X,z}$ is a normal domain, for every $z\in Z$.
\item
the natural map $\cO_{\!X}\to j_*\cO_{\!U}$ is a monomorphism,
and $\cO_{\!X}=\mathrm{i.c.}(\cO_{\!X},j_*\cO_{\!U})$.
\item
The pair $(X,Z)$ is normal.
\end{enumerate}
Then (a)$\Rightarrow$(b), since
$(j_*\cO_{\!U})_z=\cO_{\!X(z)}(U\cap X(z))$ for every $z\in Z$. Also,
(b)$\Rightarrow$(c), by virtue of \cite[Lemma 8.2.28]{Ga-Ra}.

(ii)\ \
Let $\cA$ be a quasi-coherent $\cO^a_{\!U}$-algebra. Then
$j_*\cA_{!!}$ is a quasi-coherent $\cO_{\!X}$-algebra
(\cite[Ch.I, Prop.9.4.2(i)]{EGAI} and lemma
\ref{lem_about-integral-closure}), hence the integral
closure $\cA^\nu$ of the image of $\cO^a_{\!X}$ in
$j_*\cA$ is a well defined quasi-coherent
$\cO^a_{\!X}$-algebra (see definition \ref{def_alm-int-closure}).
We call $\cA^\nu$ the {\em normalization\/} of $\cA$
over $X$.

(iii)\ \
Let $(X,Z)$ be any normal pair, and $U\subset X$ any
open subset. Directly from the definition, we see that
$(U,Z\cap U)$ is still a normal pair.
\end{remark}

\begin{lemma}\label{lem_almost-pure-crit}
Let $(X,Z)$ be a normal pair as in definition
{\em\ref{def_al-pure}(ii)}, and set $U:=X\!\setminus\!Z$.
Let $\cA$ be an \'etale $\cO^a_{\!U}$-algebra whose underlying
$\cO^a_{\!U}$-module is almost finitely presented  (resp. is
of finite rank). The following conditions are equivalent :
\begin{enumerate}
\alphaenu
\item
The normalization $\cA^\nu$ of $\cA$ over $X$ (see remark
{\em\ref{rem_almost-pairs}(ii)}) is a weakly unramified
$\cO^a_{\!X}$-algebra.
\item
$\cA^\nu$ is an \'etale $\cO^a_{\!X}$-algebra, and an
almost finitely presented $\cO^a_{\!X}$-module (resp.
and an $\cO^a_{\!X}$-module of finite rank).
\end{enumerate}
\end{lemma}
\begin{proof} Obviously (b)$\Rightarrow$(a), hence suppose
that (a) holds, and let $\cB:=\cA^\nu\otimes_{\cO^a_{\!X}}\cA^\nu$.
We may assume that $X$ is affine, and we set
$$
A:=\cA^\nu(X)
\qquad
B:=\cB(X).
$$
By assumption, the multiplication morphism $\mu:\cB\to\cA^\nu$
is flat, especially $A$ is a flat $B$-module. Let also $B'$ be
the image of $B$ in $\cB(U)$; clearly $\mu(X):B\to A$ factors
through an epimorphism $\mu':B'\to A$, therefore $A=A\otimes_BB'$,
so $A$ is also a flat $B'$-module. Pick a finite covering
$U=U_1\cup\cdots\cup U_k$ consisting of affine open subsets
of $U$. The induced morphism
$$
B'\to\cB(U_1)\times\cdots\times\cB(U_k)
$$
is a monomorphism. On the other hand, notice that
$A\otimes_{B'}\cB(U_i)=\cA(U_i)$ is an almost finitely presented
$\cB(U_i)$-module. It follows that $A$ is an almost finitely
generated projective $B'$-module
(\cite[Prop. 2.4.18 and 2.4.19]{Ga-Ra}), and therefore the kernel
of $\mu'$ is generated by an idempotent $e\in B'_*$
(\cite[Rem. 3.1.8]{Ga-Ra}). A simple inspection shows that
$e$ is necessarily the diagonal idempotent of the unramified
$\cO^a_{\!U}$-algebra $\cA$. Then all the assumptions of
proposition \ref{prop_axiomatize-familiar} are fulfilled,
so $\cA^\nu$ is an \'etale $\cO^a_{\!X}$-algebra, and an
almost finitely presented $\cO^a_{\!X}$-module.

Lastly, suppose that $\cA$ is an $\cO^a_{\!U}$-module of
finite rank, and pick $r\in\N$ such that, for every $i=1,\dots,k$,
the $r$-th exterior power of the $\cO^a_{\!U}(U_i)$-module
$\cA(U)=A\otimes_{\cO^a_{\!X}(X)}\cO^a_{\!U}(U_i)$ vanishes.
Since the induced morphism
$\cO^a_{\!X}\to\cO^a_{\!U}(U_1)\times\cdots\times\cO^a_{\!U}(U_1)$
is a monomorphism, and $A$ is a flat $\cO^a_{\!X}$-module,
we deduce that the $r$-th exterior power of $A$ vanishes
as well. Especially, $\cA^\nu$ is an $\cO^a_{\!X}$-module
of finite rank.
\end{proof}

\begin{lemma}\label{lem_pure-almost-crit}
Let $(X,Z)$ be a normal pair, $j:X\!\setminus\!Z\to X$ the
open immersion, $\cB$ any \'etale almost finitely presented
$\cO^a_{\!X}$-algebra, and $Z'\subset Z$ any constructible
closed subset. We have :
\begin{enumerate}
\item
The natural map $\cB\to j_*j^*\cB$ factors through an
isomorphism $\cB\isom(j^*\cB)^\nu$.
\item
The restriction functor \eqref{eq_resr_etale} is fully
faithful.
\item
The pair $(X,Z')$ is normal.
\item
If the pair $(X,Z)$ is almost pure, the same holds for $(X,Z')$.
\end{enumerate}
\end{lemma}
\begin{proof}(i): Set $U:=X\!\setminus\!Z$ and
$\cA:=\cB_{|U}$; the natural map $\cB\to j_*\cA$ factors
through a morphism $\phi:\cB\to\cA^\nu$ of $\cO_{\!X}^a$-algebras
(lemma \ref{lem_about-integral-closure} and remark
\ref{rem_alm-int-closure}(ii)). Fix $z\in Z$, and set
$R:=\cO^a_{\!X,z}$, $R':=(j_*\cO^a_{\!X})_z$. It follows
that the stalk $\cB_z$ is an \'etale $R$-algebra. Notice
that $(j_*\cA)_z=R'\otimes_R\cB_z$, therefore $\phi$ is
a monomorphism. Moreover, since $\cB$ is also an almost
finitely presented $\cO^a_{\!X}$-module, we have
$$
\cB_z=\mathrm{i.c.}(\cB_z,R'\otimes_R\cB_z)
\qquad
\cA^\nu_z=\mathrm{i.c.}(R,R'\otimes_R\cB_z)
$$
(proposition \ref{prop_8231}). On the other hand, $\cB_z$
is an integral $R$-algebra, therefore $\cB_z=\cA^\nu_z$,
and since $z$ is arbitrary, we conclude that $\phi$ is
an isomorphism, as asserted.

(ii): Let $\cB_1$ and $\cB_2$ be two \'etale almost finitely
presented $\cO^a_{\!X}$-algebra, set $\cA_i:=\cB_{_i|U}$ for
$i=1,2$, and let $\psi:\cA_1\to\cA_2$ be any morphism of
$\cO^a_{\!U}$-algebras; by (i), $\psi$ extends uniquely to
a morphism $\psi^\nu:\cB_1=\cA_1^\nu\to\cA^\nu_2=\cB_2$ of
$\cO^a_{\!X}$-algebras, whence the contention.

(iii): Set $U':=X\!\setminus\!Z'$ and let $j':U'\to X$ be the
open immersion; the natural morphism $\cO^a_{\!X}\to j_*\cO_{\!U}$
factors through the morphism $\cO^a_{\!X}\to j'_*\cO_{\!U'}$,
so the latter is a monomorphism. In order to show that
$(X,Z')$ is normal, it then suffices to check that the image
of $\cO^a_{\!X}$ is integrally closed in $j'_*\cO^a_{\!U'}$,
and since $\cO^a_{\!X}$ is integrally closed in $j_*\cO^a_{\!U}$,
we are reduced to proving that the natural morphism
$j'_*\cO^a_{\!U'}\to j_*\cO^a_{\!U}$ is a monomorphism.
However, let $V\subset X$ be any affine open subset,
and write $U'\cap V=V_1\cap\cdots\cap V_n$ for certain
affine open subsets $V_1,\dots,V_n$ of $X$; by assumption,
the restriction map $\cO^a_{\!X}(V_i)\to\cO^a_{\!X}(U\cap V_i)$
is a monomorphism for every $i=1,\dots,n$, so the same
holds for the restriction map
$\cO^a_{\!X}(U'\cap V)\to\cO^a_{\!X}(U\cap V)$, whence
the assertion.

(iv): Suppose that $(X,Z)$ is almost pure; in view of (ii),
in order to check that $(X,Z')$ is almost pure, it suffices
to show that every \'etale $\cO^a_{\!U'}$-algebra $\cA$ of
finite rank extends to an \'etale $\cO^a_{\!X}$-algebra of
finite rank. However, the assumption implies that
$\cA_{|U}$ extends to an $\cO^a_{\!X}$-algebra $\cB$ as
sought, and since $(X\!\setminus\!Z',Z\!\setminus\!Z')$
is normal (remark \ref{rem_almost-pairs}(iii)), (ii) says
that the isomorphism $\cB_{|U}\isom\cA$ extends to an
isomorphism $\cB_{|U'}\isom\cA$, {\em i.e.} $\cB$ is an
extension of $\cA$, as required.
\end{proof}

\begin{lemma}\label{lem_product-tail}
Let $(A_i~|~i\in\N)$ be a system of $V^a$-algebras, and
set $A:=\prod_{i\in\N}A_i$. Let also $P$ be an $A$-module,
$B$ an $A$-algebra, and suppose that
\begin{enumerate}
\alphaenu
\item
$\lim_{i\to\infty}\Ann_{V^a}(A_i)=V^a$ for the uniform
structure of \cite[Def.2.3.1]{Ga-Ra}.
\item
For every $i\in\N$, the $A_i$-modules $P_i:=P\otimes_AA_i$,
$B_i:=B\otimes_AA_i$ are almost projective of finite constant
rank equal to $i$, and $B_i$ is an \'etale $A_i$-algebra.
\end{enumerate}
Then $P$ is an almost projective $A$-module of almost finite
rank, and $B$ is an \'etale $A$-algebra.
\end{lemma}
\begin{proof} For every $j\in\N$, the finite product
$P_{\leq j}:=\prod_{i=1}^jP_i$ is an almost projective
$A$-module of finite rank, and clearly the induced morphism
$\pi_j:P\to P_{\leq j}$ is an epimorphism.
On the other hand, the $(j+1)$-th exterior power of $P$
equals the $(j+1)$-th exterior power of $\Ker\,\pi_j$,
and from condition (i) we see that
$\lim_{j\to\infty}\Ann_{V^a}(\Ker\,\pi_j)=0$, whence the
assertion for $P$.

It follows already that $B$ is an almost projective $A$-module.
It remains to show that $B$ is an unramified $A$-algebra, to
which aim, we may apply the criterion of \cite[Prop.3.1.4]{Ga-Ra}.

\begin{claim}\label{cl_haveityourway}
Under the assumptions of the lemma, the natural morphism
$$
\phi:B\otimes_AB\to C:=\prod_{i\in\N}B_i\otimes_{A_i}B_i
$$
is an isomorphism of $A$-algebras.
\end{claim}
\begin{pfclaim} For every $j\in\N$, let
$\pi_j:C\to\prod_{i\leq j} B_i\otimes_AB_i$ be the natural
morphism. Then
$$
\lim_{j\to\infty}\Ann_{V^a}\Ker(\pi_j\circ\phi)=V^a=
\lim_{j\to\infty}\Ann_{V^a}\Ker\,\pi_j.
$$
The first identity implies that $\phi$ is a monomorphism.
Next, since $\pi_j\circ\phi$ is an epimorphism, the natural
morphism $\Ker\,\pi_j\to\Coker\,\phi$ is an epimorphism;
then the second identity implies that $\phi$ is also
an epimorphism.
\end{pfclaim}

Now, by \cite[Prop.3.1.4]{Ga-Ra}, for every $i\in\N$ there
exists an idempotent $e_i\in(B_i\otimes_{A_i}B_i)_*$ uniquely
characterized by the conditions (i)--(iii) of {\em loc.cit.}
In view of claim \ref{cl_haveityourway}, the sequence
$(e_i~|~i\in\N)$ defines an idempotent in $(B\otimes_AB)_*$,
which clearly fulfills the same conditions, whence the
contention, again by \cite[Prop.3.1.4]{Ga-Ra}.
\end{proof}

\begin{proposition}\label{prop_almost-pure-crit}
Let $(X,Z)$ be a normal pair, and set $U:=X\!\setminus\!Z$.
Then the following conditions are equivalent :
\begin{enumerate}
\alphaenu
\item
The pair $(X,Z)$ is almost pure.
\item
The restriction functor
$$
\cO^a_{\!X}\Et_\mathrm{fr}\to\cO^a_{\!U}\Et_\mathrm{fr}
\quad :\quad \cA\mapsto\cA_{|U}
$$
from the category of \'etale $\cO^a_{\!X}$-algebras of finite
rank, to the category of \'etale $\cO^a_{\!U}$-algebras
of finite rank, is an equivalence.
\item
For every \'etale $\cO_{\!U}^a$-algebra $\cA$ of finite rank,
the normalization $\cA^\nu$ of $\cA$ over $X$ is an \'etale
$\cO_{\!X}^a$-algebra of finite rank (see remark
{\em\ref{rem_almost-pairs}(ii)}).
\item
For every \'etale $\cO_{\!U}^a$-algebra $\cA$ of finite rank,
$\cA^\nu$ is a weakly unramified $\cO_{\!X}^a$-algebra.
\item
Every \'etale $\cO_{\!U}^a$-algebra $\cA$ of finite rank extends
to an \'etale almost finite $\cO_{\!X}^a$-algebra.
\end{enumerate}
\end{proposition}
\begin{proof} 
The equivalence (b)$\Leftrightarrow$(c) is already
clear from lemma \ref{lem_pure-almost-crit}(i,ii).

Next, clearly (c)$\Rightarrow$(d). Conversely,
suppose that (d) holds, and let $\cA$ be any \'etale
$\cO^a_{\!U}$-algebra of finite rank, so $\cA^\nu$ is a
weakly unramified $\cO_{\!X}^a$-algebra. Then $\cA^\nu$
is actually an \'etale $\cO^a_{\!X}$-algebra, and an
almost finitely presented $\cO^a_{\!X}$-module (lemma
\ref{lem_almost-pure-crit}). Fix any
affine open subset $V\subset X$, and pick a finite covering
of $U\cap V$ consisting of affine open subsets $V_1,\dots,V_k$
of $U$. Let $r\in\N$ be an integer such that the $r$-th exterior
power of the $\cO_{\!X}^a(V_i)$-modules $\cA(V_i)$ vanish.
The induced map
$$
B:=\cO_{\!X}^a(V)\to B':=
\cO_{\!X}^a(V_1)\times\cdots\times\cO_{\!X}^a(V_k)
$$
is a monomorphism, since $(X,Z)$ is a normal pair. Moreover,
the $r$-th exterior power of the $B'$-module
$\cA^\nu(V)\otimes_BB'$ vanishes; since $\cA^\nu(V)$ is a
flat $B$-module, we conclude that the $r$-th exterior power
of $\cA^\nu(V)$ vanishes as well. Especially, $\cA^\nu$ is
of finite rank, so (c) holds.

Since (c)$\Rightarrow$(e), we suppose that (e) holds,
and deduce that (d) holds as well. Indeed, let $\cB$ be an
almost finite \'etale $\cO^a_{\!X}$-algebra extending $\cA$. In
the foregoing, we have already remarked that $\cB\subset\cA^\nu$;
especially, the diagonal idempotent of $\cA$ lies in the
image of the restriction map
$\cA^\nu\otimes_{\cO^a_{\!X}}\cA^\nu(X)_*\to
\cA\otimes_{\cO^a_{\!X}}\cA(U)_*$. Then proposition
\ref{prop_axiomatize-familiar} implies that $\cA^\nu$ is an
\'etale almost finite $\cO^a_{\!X}$-algebra, as sought.

Obviously, (a)$\Rightarrow$(e); to conclude, it suffices
then to check that (c)$\Rightarrow$(a). Thus, let $\cA$ be
an \'etale $\cO_{\!U}^a$-algebra of almost finite rank,
$V\subset U$ an affine open subset, and set $A:=\cO^a_{\!X}(V)$,
$B:=\cA(V)$. According to \cite[Th.4.3.28]{Ga-Ra}, there
exists a decomposition of $A$ as an infinite product of
a system of $V^a$-algebras $(A_i~|~i\in\N)$ fulfilling
condition (a) of lemma \ref{lem_product-tail}.
Such a decomposition determines a system of idempotent
elements $e_{V,i}\in A_*$, for every $i\in\N$, such that
$e_{V,i}\cdot e_{V,j}=0$ whenever $i\neq j$, and
characterized by the identities $e_{V,i}A=A_i$ for every $i\in\N$.
Since $\cO_{\!X*}$ is a sheaf (\cite[\S5.5.4]{Ga-Ra}),
condition (ii) ensures that, for any fixed $i\in\N$, and
$V$ ranging over the affine open subsets of $U$, the
idempotents $e_{V,i}$ glue to a global section
$e_i\in\cO^a_{\!X}(U)_*$, which will still be an idempotent;
moreover, the direct factor $\cA_i:=e_i\cA$ of $\cA$ is an
\'etale $\cO_{\!U}^a$-algebra of finite rank, and the projection
$\cA\to\cA_i$ is a morphism of $\cO^a_{\!U}$-algebras.
Especially, $A$ and $B$ fulfill as well condition (b) of
lemma \ref{lem_product-tail}.
By (c), we deduce that the normalization $\cA^\nu_i$ of
$\cA_i$ over $X$ is an \'etale $\cO^a_{\!X}$-algebra of
finite rank, for every $i\in\N$. On the other hand, notice
that the natural morphism
$$
\cO^a_{\!X*}\to j_*(\cO^a_{\!U*})=(j_*\cO^a_{\!U})_*
$$
is a monomorphism, and its image is integrally closed
in $j_*\cO^a_{\!U*}$ (\cite[Rem.8.2.30(i)]{Ga-Ra}).
It follows easily that $e_i\in\cO^a_{\!X*}(X)$;
hence
$$
e_i\cA^\nu_i=\cA^\nu_i
\qquad
\text{for every $i\in\N$}
$$
and the latter is the integral closure of $e_i\cO_{\!X}^a$
in $j_*(\cA_i)$. Consequently, $\cA^\nu=\prod_{i\in\N}\cA^\nu$
(\cite[Rem.8.2.30(ii)]{Ga-Ra}).
By lemma \ref{lem_product-tail}, it follows that $\cA^\nu$
is an \'etale $\cO^a_{\!X}$-algebra of almost finite rank.
This shows that the functor \eqref{eq_resr_etale} is
essentially surjective; combining with lemma
\ref{lem_pure-almost-crit}(ii), we get (a).
\end{proof}

\begin{corollary}\label{cor_shrink-locus}
Let $(X,Z)$ be a normal pair, and suppose that :
\begin{enumerate}
\alphaenu
\item
The scheme $Z(z)$ is finite dimensional, for every $z\in Z$.
\item
The pair $(X(z),\{z\})$ is almost pure, for every $z\in Z$.
\end{enumerate}
Then the pair $(X,Z)$ is almost pure.
\end{corollary}
\begin{proof} Let $\cA$ be any \'etale $\cO^a_{\!U}$-algebra
of finite rank, and $\cA^\nu$ the normalization of $\cA$
over $X$. In light of proposition \ref{prop_almost-pure-crit},
it suffices to show that $\cA^\nu_z$ is a weakly unramified
$\cO^a_{\!X,z}$-algebra, for every $z\in Z$. Suppose, by
way of contradiction, that the latter assertion fails;
in view of condition (a), we may then find a point $z\in Z$,
such that $\cA^\nu_z$ is not weakly unramified over
$\cO^a_{\!X,z}$, but for every proper generization $w\in Z$
of $z$, the $\cO^a_{\!X,w}$-algebra $\cA^\nu_w$ is weakly
unramified.

Now, set $U(z):=U\cap X(z)$, $V(z):=X(z)\!\setminus\!\{z\}$,
and let $j:V(z)\to U$ be the induced morphism. Notice that
$(V(z),Z(z))$ is a normal pair, and our assumption implies
that $j^*\cA^\nu$ is a weakly unramified
$\cO^a_{\!V(z)}$-algebra. By lemma \ref{lem_almost-pure-crit},
it follows that $j^*\cA^\nu$ is actually an \'etale
$\cO^a_{\!V(z)}$-algebra of finite rank. Since by assumption,
$(X(z),\{z\})$ is almost pure, proposition
\ref{prop_almost-pure-crit} says that $\cA^\nu_z$ is an
\'etale $\cO^a_{\!X,z}$-module of finite rank, a contradiction.
\end{proof}

\begin{proposition}\label{prop_pro-smooth-desc-pure}
Let $(X,Z)$ be a normal pair, and $f:X'\to X$ a morphism
of $S$-schemes; set $Z':=f^{-1}Z$, and suppose that :
\begin{enumerate}
\alphaenu
\item
$f(Z')=Z$ and $f$ is pro-smooth at every point of\/ $Z'$
(see definition {\em\ref{def_pro-smooth}}).
\item
The pair $(X',Z')$ is almost pure.
\end{enumerate}
Then the pair $(X,Z)$ is almost pure.
\end{proposition}
\begin{proof} Set $U:=X\!\setminus\!Z$, $U':=X'\!\setminus\!Z'$,
and let $f_{|U}:U'\to U$ be the restriction of $f$. Notice that
the open immersion $j':U'\to X'$ is quasi-compact
(\cite[Ch.I, Prop.6.6.4]{EGAI}), and since $f$ is flat
at every point of $Z'$ (corollary \ref{cor_pro-smooth}(i)),
$U'$ is dense in $X'$. Moreover :

\begin{claim}\label{cl_che-fatica}
(i)\ \ The pair $(X',Z')$ is normal.
\begin{enumerate}
\addenu
\item
More generally, let $\cA$ be any quasi-coherent
$\cO^a_{\!U}$-algebra. Then the natural morphism
$$
f^*(\cA^\nu)\to(f_{|U}^*\cA)^\nu
$$
is an isomorphism (notation of remark \ref{rem_almost-pairs}(ii)).
\end{enumerate}
\end{claim}
\begin{pfclaim} Since $f$ is flat at every point of $Z'$,
and the unit of adjunction $\cO^a_{\!X}\to j_*\cO^a_{\!U}$
is a monomorphism, the induced morphism
$\cO^a_{\!X'}\to f^*j_*\cO^a_{\!U}=j'_*\cO^a_{\!U'}$ is a monomorphism
as well (corollary \ref{cor_base-change-where}). Also,
$(\cO^a_{\!U})^\nu=\cO^a_{\!X}$, since $(X,Z)$ is a normal pair;
hence, (ii) implies (i).

(ii): Denote $\cA_{!!}^\nu$ the integral closure in $j_*(\cA_{!!})$
of the image of of $\cO_{\!X}$; by corollary \ref{cor_base-change-where},
the natural map
$$
f^*j_*(\cA_{!!})\to j'_*f_{|U}^*(\cA_{!!})
$$
is an isomorphism of $\cO_{\!X'}$-algebras. On the other hand,
corollary \ref{cor_pro-smooth}(ii) implies that the integral
closure of the image of $\cO_{\!X'}$ in $f^*j_*(\cA_{!!})$
equals $f^*(\cA_{!!}^\nu)$. The assertion follows.
\end{pfclaim}

Suppose now that $\cA$ is an \'etale $\cO^a_{\!U}$-algebra
of finite rank. Since, by assumption, $(X',Z')$ is almost pure,
proposition \ref{prop_almost-pure-crit} says that
$(f_{|U}^*\cA)^\nu$ is an \'etale $\cO_{\!X'}^a$-algebra
of finite rank. By claim \ref{cl_che-fatica} and lemma
\ref{lem_complem-qcoh-alg}(i) and corollary \ref{cor_pro-smooth}(ii),
we deduce that $\cA^\nu$ is a weakly \'etale $\cO_{\!X}^a$-algebra.
To conclude the proof, it suffices now to invoke proposition
\ref{prop_almost-pure-crit}.
\end{proof}

\begin{lemma}\label{lem_tight-covers}
Let $(A,I)$ be a tight henselian pair (see
\cite[Def.5.1.9]{Ga-Ra}). Then we have :
\begin{enumerate}
\item
The base change functor
$$
\bCov(\Spec\,A)\to\bCov(\Spec\,A/I)
$$
is an equivalence (notation of\/ \cite[\S8.2.22]{Ga-Ra}).
\item
More generally, the base change functor $A\Et\to A/I\Et$
restricts to an equivalence
$$
A\Et_\mathrm{afr}\isom A/I\Et_\mathrm{afr}
$$
on the respective full subcategories of \'etale algebras
of almost finite rank.
\end{enumerate}
\end{lemma}
\begin{proof} In view of \cite[Th.5.5.7(iii)]{Ga-Ra}, the
assertions follow  from lemma \ref{lem_detect-afr-after-quot}(ii).
\end{proof}

\sset\subsubsection{}\label{subsec_tight-hensel}
Let $R$ be a $V$-algebra, $I\subset R$ a principal ideal
generated by a regular element, $R^\wedge$ the $I$-adic
completion of $R$. Set
$$
X:=\Spec\,R
\qquad
X^\wedge:=\Spec\,R^\wedge
\qquad
Z:=\Spec\,R/I
\qquad
Z^\wedge:=\Spec\,R^\wedge/IR^\wedge.
$$

\begin{proposition}\label{prop_pass-to-complet}
In the situation of \eqref{subsec_tight-hensel}, suppose
furthermore that :
\begin{enumerate}
\item
There exist $n\in\N$ and a finitely generated subideal
$\fm_0\subset\fm$ such that $I^n\subset\fm_0R$.
\item
The pair $(R,I)$ is henselian.
\end{enumerate}
Then the pair $(X,Z)$ is almost pure if and only if the same
holds for the pair $(X^\wedge,Z^\wedge)$.
\end{proposition}
\begin{proof} Under the current assumptions, the pair $(R^a,I^a)$
is tight henselian (\cite[\S 5.1.12]{Ga-Ra}). Hence the assertion
is a straightforward consequence of \cite[Prop.5.4.53]{Ga-Ra}
and lemma \ref{lem_tight-covers}(i) (details left to the
reader).
\end{proof}

\sset\subsubsection{}\label{subsec_with-cotangent}
Let $p\in\Z$ be a prime integer with $pV\subset\fm$, and such
that $p$ is a regular element of $V$, and suppose that $V$
contains an element, denoted $p^{1/p}$, whose $p$-th power
generates the ideal $pV$. Let $A$ and $B$ be two flat
$V$-algebras, and for every $n\in\N$ set $V_n:=V/p^{n+1}V$,
$A_n:=A/p^{n+1}A$, $B_n:=B/p^{n+1}B$. Let also
\set\begin{equation}\label{eq_to-be-lifted}
A_0\isom B_0
\end{equation}
be a given isomorphism of $V_0$-algebras, and suppose that the
Frobenius endomorphisms of $V_0$ and $A_0$ induce isomorphisms
\set\begin{equation}\label{eq_Frob-induced}
V/p^{1/p}V\isom V_0
\qquad
A/p^{1/p}A\isom A_0.
\end{equation}
Set $X:=\Spec\,A$, $X_0:=\Spec\,A_0$, $Y:=\Spec\,B$ and
$Y_0:=\Spec\,B_0$.

\begin{corollary} In the situation of \eqref{subsec_with-cotangent},
suppose moreover that the pairs $(A,pA)$ and $(B,pB)$ are henselian.
Then the pair $(X,X_0)$ is almost pure if and only if the same
holds for the pair $(Y,Y_0)$.
\end{corollary}
\begin{proof} In view of proposition \ref{prop_pass-to-complet},
it suffices to show that \eqref{eq_to-be-lifted} lifts to
an isomorphism $A^\wedge\isom B^\wedge$ between the $p$-adic
completions of $A$ and $B$. In turns, this reduces to exhibiting
a system of isomorphisms $(\phi_n:A_n\isom B_n~|~n\in\N)$
such that $\phi_n\otimes_{V_n}V_{n_1}=\phi_{n-1}$ for every $n>0$.
To this aim, it suffices to show that
$$
\L_{A_n/V_n}=0
\qquad
\text{in $\sD(s.A\Mod)$ for every $n\in\N$}
$$
(\cite[Prop.3.2.16]{Ga-Ra}). However, in view of
\cite[Th.2.5.36]{Ga-Ra}, for every $n\in\N$, the short
exact sequence $0\to p^{n+1}V/p^{n+2}V\to V_{n+1}\to V_n\to 0$
induces a distinguished triangle in $\sD(s.A\Mod)$
$$
\L_{A_0/V_0}\to\L_{A_{n+1}/V_{n+1}}\to\L_{A_n/V_n}\to
\sigma\L_{A_0/V_0}.
$$
Hence, an easy induction further reduces to checking that
$\L_{A_0/V_0}=0$ in $\sD(s.A\Mod)$. According to
\cite[Lemma 6.5.13]{Ga-Ra}, this will follow, once we have
shown that the natural map
$$
V_{0,(\Phi)}\derotimes_{V_0}A_0\to A_{0,(\Phi)}
$$
is an isomorphism in $\sD(V_0\Mod)$ (notation of {\em loc.cit.}).
Taking into account the isomorphisms \eqref{eq_Frob-induced},
this holds if and only if the natural map
$$
V/p^{1/p}V\derotimes_{V_0}A_0\to A/p^{1/p}A
$$
is an isomorphism in $\sD(V_0\Mod)$. The latter assertion is
clear, since $A$ is a flat $V$-algebra.
\end{proof}

\subsection{Normalized lengths}\label{sec_norm-lengths}
Let $(V,|\cdot|)$ be any valuation ring, with value group
$\Gamma_{\!V}$; according to our general conventions, the
composition law of $\Gamma_{\!V}$ is denoted multiplicatively;
however, sometimes it is convenient to switch to an additive
notation. Hence, we adopt the notation :
$$
(\log\Gamma_{\!V},\leq)
$$
to denote the ordered group $\Gamma_{\!V}$ with additive
composition law and whose ordering is the reverse of the
original ordering of $\Gamma_{\!V}$.
The unit of $\log\Gamma_{\!V}$ shall be naturally denoted by $0$,
and we shall extend the ordering of $\log\Gamma_{\!V}$ by adding
a largest element $+\infty$, as customary. Also, we set :
$\log\Gamma_{\!V}^+:=\{\gamma\in\log\Gamma_{\!V}~|~\gamma\geq 0\}$
and $\log|0|:=+\infty$.

\sset\subsubsection{}\label{subsec_general-uniformity}
In this section, $(K,|\cdot|)$ denotes a valued field of rank
one, with value group $\Gamma$.
We let $\kappa$ be the residue field of $K^+$. As usual, we
set $S:=\Spec\,K^+$, and denote by $s$ (resp. $\eta$) the
closed (resp. generic) point of $S$. Let $\fm_K\subset K^+$
be the maximal ideal, and set $\fm:=\fm_K$ in case $\Gamma$
is not discrete, or else $\fm:=K^+$, in case $\Gamma\simeq\Z$;
in the following, whenever we refer to almost rings or almost
modules, we shall assume -- unless otherwise stated --
that the underlying almost ring theory is the one defined
by the standard setup $(K^+,\fm)$ (see \cite[\S6.1.15]{Ga-Ra}).

Let $A$ be any $K^{+a}$-algebra, and $c$ a cardinal
number; following \cite[\S2.3]{Ga-Ra}, we denote by $\cM_c(A)$
the set of isomorphism classes of $K^{+a}$-modules which
admit a set of generators of cardinality $\leq c$.
The set $\cM_c(A)$ carries a natural uniform structure (see
\cite[Def.2.3.1]{Ga-Ra}), which admits the fundamental system
of entourages
$$
(E_\gamma~|~\gamma\in\log\Gamma^+\setminus\{0\})
$$
defined as follows. For any $b\in\fm\setminus\{0\}$, we let
$E_{|b|}$ be the set of all pairs $(M,M')$ of elements
of $\cM_c(A)$ such that there exist a third $A$-module
$N$ and $A$-linear morphisms $N\to M$, $N\to M'$ whose
kernel and cokernel are annihilated by $b$.

\sset\subsubsection{}\label{subsec_aim-length}
The aim of this section is to define and study a well-behaved
notion of {\em normalized length\/} for torsion modules $M$ over
$K^+$-algebras of a fairly general type. This shall be achieved
in several steps. Let us first introduce the categories of
algebras with which we will be working.

\begin{definition}\label{def_working-algeb}
Let $V$ be any valuation ring, with maximal ideal $\fm_V$.

(i)\ \
We let $V\mAlg_0$ be the subcategory of $V\Alg$
whose objects are the local and essentially finitely
presented $V$-algebras $A$ whose maximal ideal $\fm_A$
contains $\fm_VA$. The morphisms $\phi:A\to B$ in
$K^+\mAlg_0$ are the local maps. Notice that every
morphism in $V\mAlg_0$ is an essentially finitely
presented ring homomorphism (\cite[Ch.IV, Prop.1.4.3(v)]{EGAIV}).

Recall as well, that every object of $V\mAlg_0$
is a coherent ring (see \eqref{sec_sch-val-rings}). More
generally, if $\underline A:=(A_i~|~i\in I)$ is any filtered
system of essentially finitely presented $V$-algebras with
flat transition morphisms, then the colimit of $\underline A$
is still a coherent ring (lemma \ref{lem_rflx-on-lim}(ii.a)).

(ii)\ \
We say that a local $V$-algebra $A$ is {\em measurable\/} if
it admits an ind-\'etale local map of $V$-algebras $A_0\to A$,
from some object $A_0$ of $V\mAlg_0$. The measurable
$V$-algebras form a category $V\mAlg$, whose morphisms
are the local maps of $V$-algebras. As noted above, every
measurable $V$-algebra is a coherent ring.
\end{definition}

\begin{lemma}\label{lem_lucilla}
In the situation of definition {\em\ref{def_working-algeb}},
let $(A,\fm_A)$ be any measurable $V$-algebra. Then the
following holds :
\begin{enumerate}
\item
$A/\fm_VA$ is a noetherian ring.
\item
If the valuation of $V$ has finite rank, every $\fm_A$-primary
ideal of $A$ contains a finitely generated $\fm_A$-primary subideal.
\item
For every $A$-module $M$ of finite type supported at the
closed point of $\Spec\,A$, the $A$-module $M/\fm_VM$ has
finite length.
\item
For every finitely generated ideal $I\subset A$, the $V$-algebra
$A/I$ is measurable.
\end{enumerate}
\end{lemma}
\begin{proof}(i): Let $A^\sh$ be the strict henselization
of $A$ at a geometric point localized at the closed point;
then $A^\sh/\fm_VA^\sh$ is a strict henselization of
$A/\fm_VA$ (\cite[Ch.IV, Prop.18.8.10]{EGA4}), and it is
therefore also a strict henselization of a $V/\fm_V$-algebra
of finite type. Hence $A^\sh/\fm_VA^\sh$ is noetherian,
and then the same holds for $A/\fm_VA$
(\cite[Ch.IV, Prop.18.8.8(iv)]{EGA4}).

(iii) follows easily from (i) : the details shall be left
to the reader.

(ii): If the valuation of $V$ has finite rank, there exists
an element $t_0\in\fm_V$ that generates a $\fm_V$-primary
ideal. Let $I\subset A$ be a $\fm_A$-primary ideal; then
$t_0^N\in I$, for some integer $N>0$; moreover, the image
$\bar I$ of $I$ in $A/\fm_VA$ is finitely generated, by
(i). Pick elements $t_1,\dots,t_n\in I$ whose images in
$A/\fm_VA$ form a system of generators of $\bar I$; it
follows that $t_0^N,t_1,\dots,t_n$ form a $\fm_A$-primary
ideal contained in $I$.

(iv): We may find an ind-\'etale local morphism $A_0\to A$
from an object $A_0$ of $V\mAlg_0$, and a finitely generated
ideal $I_0\subset A_0$ such that $I=I_0A$. Then $A/I_0$
is an object of $V\mAlg_0$ as well, and the induced map
$A_0/I_0\to A/I$ is ind-\'etale.
\end{proof}

\sset\subsubsection{}\label{subsec_replace-cor}
Now, suppose that $V$ is both a valuation ring and a flat,
measurable $K^+$-algebra, and denote $\fm_V$ (resp. $\kappa(V)$,
resp. $|\cdot|_V$) the maximal ideal (resp. the residue field,
resp. the valuation) of $V$. We claim that $V$ has rank one,
and the ramification index $(\Gamma_{\!\!V}:\Gamma)$ is finite.
Indeed, let us write $V$ as the colimit of a filtered system
$(V_i~|~i\in I)$ of objects of $K^+\mAlg_0$ with essentially
\'etale transition maps; it follows that each $V_i$ is a
valuation ring of rank one (claim \ref{cl_finally-named}).
Also, the transition maps induce isomorphisms on the value groups :
indeed, this is clear if $K^+$ is a discrete  valuation ring,
since in that case the same holds for the $V_i$, and the
transition maps are unramified by assumption; in the case
where $\Gamma$ is not discrete, the assertion follows from
corollary \ref{cor_integral-fp-ext}. Therefore, we are
reduced to the case where $V$ is an object of $K^+\mAlg_0$,
to which corollary \ref{cor_integral-fp-ext} applies.

Suppose first that $M$ is a finitely generated torsion
$V$-module; in this case the Fitting ideal $F_0(M)\subset V$
is well defined, and it is shown in
\cite[Lemma 6.3.1 and Rem.6.3.5]{Ga-Ra} that the map
$M\mapsto F_0(M)$ is additive on the set of isomorphism
classes of finitely generated $V$-modules, {\em i.e.}
for every short exact sequence $0\to M_1\to M_2\to M_3\to 0$
of such modules, one has :
\set\begin{equation}\label{eq_additivity}
F_0(M_2)=F_0(M_1)\cdot F_0(M_3).
\end{equation}
Next, the almost module $F_0(M)^a$ is an element of the
group of fractional ideals $\Div(V^a)$ defined in
\cite[\S6.1.16]{Ga-Ra}, and there is a natural isomorphism
\set\begin{equation}\label{eq_Div=Gamma}
\Div(V^a)\simeq\log\Gamma^\wedge_{\!\!V}
\qquad I\mapsto|I|
\end{equation}
where $\Gamma^\wedge_{\!\!V}$ is the completion of $\Gamma_{\!\!V}$
for the uniform structure deduced from the ordering : see
\cite[Lemma 6.1.19]{Ga-Ra}. Hence we may define :
$$
\lambda_V(M):=(\Gamma_{\!\!V}:\Gamma)\cdot|F_0(M)^a|\in
\log\Gamma^\wedge_{\!\!V}.
$$
In view of \eqref{eq_additivity}, we see that :
\set\begin{equation}\label{eq_ineq-f.g.case}
\lambda_V(M')\leq\lambda_V(M)\qquad \text{ whenever $M'\subset M$
are finitely generated.}
\end{equation}
More generally, if $M$ is any torsion $V$-module, we let :
\set\begin{equation}\label{eq_full-gener}
\lambda_V(M):=
\sup\,\{\lambda_V(M')~|~M'\subset M, \text{ $M'$ finitely generated}\}
\in\log\Gamma^\wedge_{\!\!V}\cup\{+\infty\}
\end{equation}
which, in view of \eqref{eq_ineq-f.g.case}, agrees with
the previous definition, in case $M$ is finitely generated.

\sset\subsubsection{}\label{subsec_unif-alm-fg}
Suppose now that the $V^a$-module $M^a$ is uniformly almost
finitely generated. Then according to \cite[Prop.2.3.23]{Ga-Ra}
one has a well defined Fitting ideal $F_0(M^a)\subset V^a$,
which agrees with $F_0(M)^a$ in case $M$ is finitely generated.

\begin{lemma}\label{lem_agreement}
In the situation of \eqref{subsec_unif-alm-fg}, we have :
\begin{enumerate}
\item
\ $\lambda_V(M)=(\Gamma_{\!\!V}:\Gamma)\cdot|F_0(M^a)|$.
\item
If\/ $0\to N_1\to N_2\to N_3\to 0$ is any short exact sequence
of uniformly almost finitely generated $V^a$-modules, then
$$
|F_0(N_2)|=|F_0(N_1)|+|F_0(N_3)|.
$$
\end{enumerate}
\end{lemma}
\begin{proof} (ii): It is a translation of
\cite[Lemma 6.3.1 and Rem.6.3.5(ii)]{Ga-Ra}.

(i): Let $k$ be a uniform bound for $M^a$, and
denote by $\cI_k(M)$ the set of all submodules of $M$
generated by at most $k$ elements. Set $e:=(\Gamma_{\!\!V}:\Gamma)$;
if $N\in\cI_k(M)$, then $\lambda_V(N)\leq e\cdot|F_0(M^a)|$,
due to (ii). By inspecting the definitions, it then follows that
$$
e\cdot|F_0(M^a)|=
\sup\,\{\lambda_V(N)~|~N\in\cI_k(M)\}\leq\lambda_V(M).
$$
Now, suppose $M'\subset M$ is any submodule generated by,
say $r$ elements; for every $\eps\in\fm$ we can find
$N\in\cI_k(M)$ such that $\eps M\subset N$, hence
$\eps M'\subset N$, therefore
$\lambda_V(\eps M')\leq\lambda_V(N)\leq e\cdot|F_0(M^a)|$.
However,
$\lambda_V(M')-\lambda_V(\eps M')=\lambda_V(M'/\eps M')\leq
\lambda_V(K^{+\oplus r}/\eps K^{+\oplus r})=
r\cdot\lambda_V(K^+/\eps K^+)$. We deduce easily that
$\lambda_V(M')\leq e\cdot|F_0(M^a)|$, whence the claim.
\end{proof}

\begin{proposition}\label{prop_additive-vr}
{\em(i)}\ \ Let\/ $0\to M_1\to M_2\to M_3\to 0$
be a short exact sequence of torsion $V$-modules. Then :
$$
\lambda_V(M_2)=\lambda_V(M_1)+\lambda_V(M_3).
$$
\begin{enumerate}
\addenu
\item
$\lambda_V(M)=0$ if and only if $M^a=0$.
\item
Let $M$ be any torsion $V$-module. Then :
\begin{enumerate}
\item
If $M$ is finitely presented, $\lambda_V(N)>0$ for every
non-zero submodule $N\subset M$.
\item
If $(M_i~|~i\in I)$ is a filtered system of submodules
of $M$, then :
$$
\lambda_V(\colim_{i\in I}M_i)=\lim_{i\in I}\,\lambda_V(M_i).
$$
\end{enumerate}
\end{enumerate}
\end{proposition}
\begin{proof} (ii): By lemma \ref{lem_agreement}(i) it is
clear that $\lambda_V(M)=0$ whenever $M^a=0$. Conversely,
suppose that $\lambda_V(M)=0$, and let $m\in M$ be any
element; then necessarily $\lambda_V(Vm)=0$, and then
it follows easily that $\fm\subset\Ann_V(m)$, so
$M^a=0$.

(iii.a): It suffices to show that $M$ does not contain
non-zero elements that are annihilated by $\fm$, which
follows straightforwardly from \cite[Lemma 6.1.14]{Ga-Ra}.

(iii.b): Set $M:=\colim_{i\in I}M_i$. By inspecting the
definitions we see easily that
$\lambda_V(M)\geq\sup\{\lambda_V(M_i)~|~i\in I\}$.
To show the converse inequality, let $N\subset M$ be any
finitely generated submodule; we may find $i\in I$ such that
$N\subset M_i$, hence $\lambda_V(N)\leq\lambda_V(M_i)$, and
the assertion follows.

(i): Set $e:=(\Gamma_{\!\!V}:\Gamma)$. We shall use the following :

\begin{claim}\label{cl_sub-of-uni}
Every submodule of a uniformly almost finitely generated
$V^a$-module is uniformly almost finitely generated.
\end{claim}
\begin{pfclaim} Let $N'\subset N$, with $N$ uniformly
almost finitely generated, and let $k$ be a uniform bound
for $N$; for every $\eps\in\fm$ we can find $N''\subset N$
such that $N''$ is generated by at most $k$ almost elements
and $\eps N\subset N''$. Clearly, it suffices to show that
$N'\cap N''$ is uniformly almost finitely generated and admits
$k$ as a uniform bound, so we may replace $N$ by $N''$
and $N'$ by $N'\cap N''$, and assume from start that
$N$ is finitely generated. Let us pick an epimorphism
$\phi:(V^a)^{\oplus k}\to N$; it suffices to show that
$\phi^{-1}(N')$ is uniformly almost finitely generated
with $k$ as a uniform bound, so we are further reduced
to the case where $N$ is free of rank $k$. Then we can write
$N'=L^a$ for some submodule $L\subset V^{\oplus k}$;
notice that $(V^{\oplus k}/L)^a$ is almost finitely presented,
since it is finitely generated (\cite[Prop.6.3.6(i)]{Ga-Ra}),
hence $N'$ is almost finitely generated (\cite[Lemma 2.3.18(iii)]{Ga-Ra}).
Furthermore, $L$ is the colimit of the family $(L_i~|~i\in I)$
of its finitely generated submodules, and each $L_i$ is
a free $V$-module (\cite[Ch.VI, \S3, n.6, Lemma 1]{BouAC}).
Necessarily the rank of $L_i$ is $\leq k$ for every $i\in I$,
hence $\Lambda^{k+1}_VL\simeq
\colim_{i\in I}\Lambda^{k+1}_VL_i=0$, and then the claim
follows from \cite[Prop.6.3.6(ii)]{Ga-Ra}.
\end{pfclaim}

Now, let $N\subset M_2$ be any finitely generated
submodule, $\bar N\subset M_3$ the image of $N$; by
claim \ref{cl_sub-of-uni}, $M_1\cap N$ is uniformly
almost finitely generated, hence lemma \ref{lem_agreement}(i,ii)
shows that :
$$
\lambda_V(N)=e\cdot|F_0(N^a)|=
e\cdot|F_0(M_1^a\cap N^a)|+e\cdot|F_0(\bar N{}^a)|=
\lambda_V(M_1\cap N)+\lambda_V(\bar N).
$$
Taking the supremum over the family $(N_i~|~i\in I)$ of all
finitely generated submodules of $M_2$ yields the identity :
$$
\lambda_V(M_2)=\lambda_V(M_3)+\sup\{\lambda_V(M_1\cap N_i)~|~i\in I\}.
$$
By definition, $\lambda_V(M_1\cap N_i)\leq\lambda_V(M_1)$ for every
$i\in I$; conversely, every finitely generated submodule of $M_1$
is of the form $N_i$ for some $i\in I$, whence the contention.
\end{proof}

\begin{remark}\label{rem_shallbeleft}
Suppose that $\Gamma\simeq\Z$, and let $\gamma_0\in\log\Gamma^+$
be the positive generator. Then by a direct inspection of the
definition one finds the identity :
$$
\lambda(M)=(\Gamma_{\!\!V}:\Gamma)\cdot\length_V(M)\cdot\gamma_0
$$
for every torsion $V$-module $M$. The verification shall be left
to the reader.
\end{remark}

\sset\subsubsection{}\label{subsec_element-divs}
Let $M$ be a torsion $V$-module, such that $M^a$ is almost
finitely generated; we wish now to explain that $\lambda_V(M)$
can also be computed in terms of a suitable sequence of elementary
divisors for $M^a$. Indeed, suppose first that $M$ is a finitely
presented torsion $V$-module; then we have a decomposition
\set\begin{equation}\label{eq_tyger}
M=(V/a_0V)\oplus\cdots\oplus(V/a_nV)
\qquad
\text{where $n:=\dim_{\kappa(V)} M/\fm_V M-1$}
\end{equation}
for certain $a_0,\dots,a_n\in\fm_V$ (\cite[Lemma 6.1.14]{Ga-Ra}).
Clearly $\gamma_i:=\log|a_i|_V>0$ for every $i=0,\dots,n$, and
after reordering we may assume that $\gamma_0\geq\cdots\geq\gamma_n$;
then we may set $\gamma_i:=0$ for every $i>n$, and the resulting 
sequence $(\gamma_i~|~i\in\N)$ of {\em elementary divisors\/}
of $M$ is independent of the chosen decomposition
\eqref{eq_tyger}, since we have more precisely :
\set\begin{equation}\label{eq_halloween}
a_iV=\Ann_V(\Lambda_V^{i+1}M)
\qquad
\text{for every $i\in\N$}.
\end{equation}
Indeed, from a decomposition \eqref{eq_tyger} with
$a_jV\subset a_{j+1}V$ for every $j<n$, we get a $V$-linear
surjection $M\to(V/a_iV)^{\oplus i+1}$, whence a surjection
$$
N_i:=\Lambda_V^{i+1}M\to\Lambda^{i+1}_V(V/a_iV)^{\oplus i+1}
\isom V/a_iV
$$
which shows that $J_i:=\Ann_VN_i\subset a_iV$. For the converse
inclusion, notice that $N_i$ is generated by the system of all
elements of the form $\omega:=x_0\wedge\cdots\wedge x_i$, where
$x_0\in M_{j(0)},\dots,x_i\in M_{j(i)}$ for some strictly
increasing map $j:\{0,\dots,i\}\to\{0,\dots,n\}$;
especially, $j(i)\geq i$, whence $a_i\omega=0$, and the
assertion follows. Moreover, a simple inspection yields
the identity
\set\begin{equation}\label{eq_div-comp-length}
\lambda_V(M)=(\Gamma_{\!V}:\Gamma)\cdot(\gamma_0+\cdots+\gamma_n).
\end{equation}
We regard the sequence $(\gamma_i~|~i\in\N)$ as an element of the
$\log\Gamma^\wedge_{\!\!V}$-normed space $\ell^\infty(\Gamma^+_{\!\!V})$
of {\em bounded\/} sequences of elements of $\log\Gamma_{\!\!V}^{\wedge+}$,
{\em i.e.} the set of all sequences
$\underline\delta:=(\delta_i~|~i\in\N)$ with
$$
\Vert\underline\delta\Vert:=\sup(\delta_i~|~i\in\N)<+\infty.
$$

\begin{lemma}\label{lem_rachmaninov}
Let $\phi:N\to N'$ be a map of finitely presented torsion
$V$-modules, and denote by $(\gamma_i~|~i\in\N)$ (resp.
$(\gamma'_i~|~i\in\N)$) the sequence of elementary divisors
of $N$ (resp. $N'$). Then :
\begin{enumerate}
\item
If $\phi$ is injective, we have\ \ $\gamma_i\leq\gamma'_i$\ \ for
every $i\in\N$.
\item
If $\phi$ is surjective, we have\ \ $\gamma_i\geq\gamma'_i$\ \ for
every $i\in\N$.
\end{enumerate}
\end{lemma}
\begin{proof} Denote by $K_V$ the field of fractions of $V$,
and set $P^\dagger:=\Hom_V(P,K_V/V)$ for every $V$-module $P$.
We have a natural bilinear pairing $P\otimes_VP^\dagger\to K_V/V$,
whence a $V$-linear map
$$
\omega_P:P\to(P^\dagger)^\dagger.
$$

\begin{claim}\label{cl_first-for-V-mods}
(i)\ \
If $\phi$ is injective, $\phi^\dagger:=\Hom_V(\phi,K_V/V)$
is surjective.
\begin{enumerate}
\addenu
\item
For every finitely presented $V$-module $P$, we have :
\begin{enumerate}
\item
The map $\omega_P$ is an isomorphism.
\item
$P$ and $P^\dagger$ are isomorphic $V$-modules.
\end{enumerate}
\end{enumerate}
\end{claim}
\begin{pfclaim} (ii) is left to the reader.
(i) follows immediately from lemma \ref{lem_coh-inject}(i).
\end{pfclaim}

By virtue of claim \ref{cl_first-for-V-mods}(ii.b), the
$V$-modules $N$ and $N^\dagger$ have the same sequences
of elementary divisors. Taking into account claim
\ref{cl_first-for-V-mods}(i), we deduce that
(ii)$\Rightarrow$(i), so it remains only to show that
(ii) holds. However, if $\phi$ is surjective, the same
holds for $\Lambda^{i+1}_V\phi$, and then the assertion
follows immediately from \eqref{eq_halloween}.
\end{proof}

\begin{proposition} Let $0\to M'\to M\to M''\to 0$ be
a short exact sequence of finitely presented torsion
$V$-modules. Then for every $i,j,t\in\N$ we have :
$$
\sum^{i+j+t}_{k=i+j}\gamma_k(M)\leq
\sum^{i+t}_{k=i}\gamma_k(M')+\sum^{j+t}_{k=j}\gamma_k(M'').
$$
\end{proposition}
\begin{proof} Let us choose a decomposition for $M$
as in \eqref{eq_tyger}, with
$a_0V\subset a_1V\subset\cdots\subset a_nV$, and set
$\bar M:=(V/b_0V)\oplus\cdots\oplus(V/b_{i+j+t}V)$, where
$b_i=a_{i+j}$ for $i=0,\dots,i+j$ and $b_i=a_i$ for
$i=i+j+1,\dots,i+j+t$. We have an obvious $V$-linear
surjection $M\to\bar M$, and we let as well $\bar M{}'$
be the image of $M'$ in $\bar M$, and
$\bar M{}'':=\bar M/\bar M{}'$. By construction, we have
$\sum^{i+j+t}_{k=i+j}\gamma_k(M)=\sum^{i+j+t}_{k=i+j}\gamma_k(\bar M)$,
and taking into account lemma \ref{lem_rachmaninov}(ii),
we are easily reduced to checking the stated inequality
for the resulting short exact sequence
$0\to\bar M{}'\to\bar M\to\bar M{}''\to 0$. Next, for
$k=0,\dots,i+j+t$ we consider the $V$-linear map
$\phi_k:V\to V/a_{i+j}V$ that maps $1$ to the class of
$a_{i+j}/b_k$; clearly $\phi_k$ factors through a $V$-linear
injection $\bar\phi_k:V/b_kV\to V/a_{i+j}V$, and the direct
sum of the latter maps is a $V$-linear injection
$\bar\phi:\bar M\to N:=(V/a_{i+j}V)^{\oplus i+j+t+1}$.
Let $N'\subset N$ be the image of $\bar M{}'$, and
set $N'':=N/N'$; we notice that
$$
\sum^{i+j+t}_{k=i+j}(\gamma_k(N)-\gamma_k(\bar M))=
|F_0(N)|-|F_0(\bar M)|=C:=
\sum_{k=1}^t(\log|a_{i+j}|_V-\log|a_{i+j+k}|_V).
$$
Taking into account proposition \ref{prop_additive-vr}(i),
we deduce that
$$
\sum^{j+t}_{k=j}(\gamma_k(N'')-\gamma_k(\bar M{}''))\leq
|F_0(N'')|-|F_0(\bar M{}'')|=C.
$$
Summing up, we are further reduced to checking the
sought inequality for the short exact sequence
$0\to N'\to N\xrightarrow{\ \pi\ }N''\to 0$. Let us
remark, quite generally :

\begin{claim}\label{cl_halloween}
Let $A$ be a local ring, $M$ a finitely generated $A$-module,
$F,F'$ two free $A$-modules of the same rank $r\in\N$, and
$\pi:F\to M$, $\pi':F'\to M$ two $A$-linear surjections.
Then there exists an $A$-linear isomorphism $\omega:F\isom F'$
such that $\pi'\circ\omega=\pi$.
\end{claim}
\begin{pfclaim} Let $\kappa$ be the residue field of $A$,
and set $m:=\dim_\kappa(M\otimes_A\kappa)$; by Nakayama's
lemma, $M$ admits a minimal system of generators
$x_1,\dots,x_m$, and we choose $f_1,\dots,f_m\in F$ with
$\pi(f_i)=x_i$ for $i=1,\dots,m$. Next, notice that
$G:=\Ker\,(\pi\otimes_A\kappa)$ is a $\kappa$-vector
space of dimension $r-m$, and pick a system of elements
$f_{m+1},\dots,f_r\in F$ whose image in $G$ is a basis;
applying again Nakayama's lemma, it is easily seen that
the system $f_1,\dots,f_r$ generates the $A$-module $F$,
and it is therefore a basis of the latter. Likewise, we
may find a basis $f'_1,\dots,f'_r$ of $F'$ such that
$\pi'(f_i)=x_i$ for $i=1,\dots,m$ and $\pi'(f_i)=0$
for $i=m+1,\dots,r$; then the unique $A$-linear map
$\omega:F\to F'$ such that $\omega(f_i)=f'_i$ for
$i=1,\dots,r$ will do.
\end{pfclaim}

Now, say that $N''=(V/c_0V)\oplus\cdots\oplus(V/c_pV)$ is a
minimal decomposition with $c_0V\subset\cdots\subset c_pV$;
clearly we must have $p\leq i+j+t$ and $a_{i+j}\in c_0V$.
Set $c_{p+1}=\cdots=c_{i+j+t}=1$, and let
$\pi'_k:V/a_{i+j+t}\to V/c_kV$ be the natural projection,
for $k=0,\dots,i+j+t$; the direct sum of the maps $\pi'_k$
is a $V$-linear surjection $\pi':N\to N''$, and by virtue
of claim \ref{cl_halloween} there exists a $V$-linear
automorphism $\omega$ of $N$ such that $\pi'\circ\omega=\pi$.
Thus, it suffices to show the stated inequality for the
short exact sequence
$0\to Q:=\omega(N')\to N\xrightarrow{\ \pi'\ }N''\to 0$.
But notice that $Q=\bigoplus_{k=0}^{i+j+t}\ker\,\pi'_k$,
and therefore
$$
\gamma_k(Q)=\log|a_{i+j}|_V-\log|c_{i+j+t-k}|_V
\qquad
\text{for $k=0,\dots,i+j+t$}.
$$
On the other hand $\gamma_k(N)=\log|a_{i+j}|_V$ and
$\gamma_k(N'')=\log|c_k|_V$ for $k=0,\dots,i+j+t$;
summing up, we conclude the the sought inequality
holds, and indeed it is an equality for this last
short exact sequence.
\end{proof}

\sset\subsubsection{}\label{subsec_al-elem-divs}
We fix now a large cardinal number $\omega$, and write
just $\cM(V^a)$ instead of $\cM_\omega(V^a)$. Also,
for any $\gamma\in\log\Gamma^+_{\!\!V}$, set $[-\gamma,\gamma]:=
\{\delta\in\log\Gamma_{\!\!V}~|~-\gamma\leq\delta\leq\gamma\}$.
Suppose that $N$ and $N'$ are two finitely presented
$V$-modules such that $(N^a,N^{\prime a})\in E_\delta$
for some $\delta\in\log\Gamma_{\!\!V}^+\setminus\{0\}$. Say
that $\delta=\log|b|_V$ for some $b\in\fm_V$; by standard
arguments, we obtain maps $\phi:N\to N'$ and
$\phi':N'\to N$ such that $\phi'\circ\phi=b^4\cdot\one_N$
and $\phi'\circ\phi=b^4\cdot\one_{N'}$. Let now
$(\gamma_i~|~i\in\N)$ (resp. $(\gamma'_i~|~i\in\N)$) be
the sequence of elementary divisors for $N$ (resp. for $N'$).
Then the sequence of elementary divisors for $b^4N$ is 
$(\max(0,\gamma_i-4\delta)~|~i\in\N)$, and likewise for
$b^4N'$. In view of lemma \ref{lem_rachmaninov}, we deduce
easily that
$$
\gamma_i-\gamma'_i\in[-4\delta,4\delta]
\qquad
\text{for every $i\in\N$}.
$$
Consider a torsion almost finitely generated $V^a$-module $M$,
and recall that $M$ is almost finitely presented
(\cite[Prop.6.3.6(i)]{Ga-Ra}); we then may attach to $M$ a net
of elements of $\ell^\infty(\log\Gamma^+_{\!\!V})$, as follows. For
every $\delta\in\log\Gamma^+_{\!\!V}\setminus\{0\}$, pick a finitely
presented $V$-module $N_\delta$ such that
$(M,N^a_\delta)\in E_\delta$, and denote by $\underline\gamma_\delta$
the sequence of elementary divisors of $N_\delta$. The foregoing
easily implies that the system
$(\underline\gamma_\delta~|~\delta\in\log\Gamma^+_{\!\!V}\setminus\{0\})$
is a net for the uniform structure of $\ell^\infty(\log\Gamma^+_{\!\!V})$
induced by the norm $\Vert\cdot\Vert$. However,
$(\ell^\infty(\log\Gamma^+_{\!\!V}),\Vert\cdot\Vert)$ is a complete
normed space, hence this net converges to a well defined
sequence $\underline\gamma_M:=
(\gamma_i~|~i\in\N)\in\ell^\infty(\log\Gamma^+_{\!\!V})$.
It is easily seen that $\underline\gamma_M$ is independent
of the chosen net, and defines an invariant which we call
the {\em sequence of elementary divisors\/} of $M$. A simple
inspection of the construction shows that the sequence
$\underline\gamma_M$ is monotonically decreasing, and
$$
\lim_{i\to+\infty}\gamma_i=0.
$$

\sset\subsubsection{}
Denote again by $K_V$ the field of fractions of $V$; as in
the proof of lemma \ref{lem_rachmaninov}, we consider the
functor
$$
V^a\Mod\to V^a\Mod
\quad :\quad
M\mapsto M^\dagger:=\Alhom_{V^a}(M,K_V^a/V^a).
$$
For every $V^a$-module $M$ we have a natural $V^a$-bilinear
pairing
$$
M\otimes_{V^a}M^\dagger\to K_V^a/V^a
$$
which in turns yields a natural transformation
\set\begin{equation}\label{eq_double-dual}
M\to(M^\dagger)^\dagger.
\end{equation}
Together with the constructions of \eqref{subsec_al-elem-divs},
we may now extend lemma \ref{lem_rachmaninov} {\em verbatim}
to arbitrary almost finitely generated $V^a$-modules.

\begin{proposition}\label{prop_rachmaninov}
Let $\phi:M\to M'$ be a morphism of almost finitely generated
torsion $V^a$-modules, and denote by $(\gamma_i~|~i\in\N)$ (resp.
$(\gamma'_i~|~i\in\N)$) the sequence of elementary divisors
of $M$ (resp. $M'$). Then :
\begin{enumerate}
\item
$M^\dagger$ is an almost finitely generated torsion $V^a$-module,
the morphism \eqref{eq_double-dual} is an isomorphism, and the
sequences of elementary divisors of $M$ and $M^\dagger$ coincide.
\item
If\/ $\phi$ is a monomorphism, we have\ \ $\gamma_i\leq\gamma'_i$\ \
for every $i\in\N$.
\item
If\/ $\phi$ is an epimorphism, we have\ \ $\gamma_i\geq\gamma'_i$\ \
for every $i\in\N$.
\end{enumerate}
\end{proposition}
\begin{proof}(i): Since $M$ is almost finitely presented
(\cite[Prop.6.3.6(i)]{Ga-Ra}), the assertion follows
immediately from \cite[Lemma 2.3.7(iii)]{Ga-Ra} and claim
\ref{cl_first-for-V-mods}(ii).

(iii): To begin with, we remark more generally :

\begin{claim}\label{cl_directly-from-mouth}
Let $P$, $Q$ be two torsion almost finitely generated
$V^a$-modules, and $\delta\in\log\Gamma^+_{\!\!V}$ any
element such that $(P,Q)\in E_\delta$. Denote by
$(\gamma^P_i~|~i\in\N)$ (resp. $(\gamma^Q_i~|~i\in\N)$)
the sequence of elementary divisors of $P$ (resp. of $Q$).
Then
$$
\gamma^P_i-\gamma^Q_i\in[-4\delta,4\delta]
\qquad
\text{for every $i\in\N$}.
$$
\end{claim}
\begin{pfclaim} For any $\eps\in\log\Gamma^+_{\!\!V}$
pick a torsion finitely presented $V$-module $N$ such
that $(P,N^a)\in E_\eps$, and let $(\gamma^N_i~|~i\in\N)$
be the sequence of elementary divisors of $N$. It follows
that $(Q,N^a)\in E_{\delta+\eps}$, and from the construction
of \eqref{subsec_al-elem-divs} we see that
$$
\gamma^N_i-\gamma^P_i\in[-4\eps,4\eps]
\qquad
\gamma^N_i-\gamma^Q_i\in[-4(\eps+\delta),4(\eps+\delta)]
\qquad
\text{for every $i\in\N$}
$$
so $\gamma^P_i-\gamma_i^Q\in[-8\eps-4\delta,8\eps+4\delta]$
for every $i\in\N$. Since $\eps$ is arbitrary, the claim follows.
\end{pfclaim}

Fix $\delta\in\log\Gamma_{\!\!V}^+\setminus\{0\}$, pick
finitely presented $V$-modules $N$, $N'$ such that
$$
(N^a,M),(N'{}^a,M')\in E_\delta
$$
and denote by $(\beta_i~|~i\in\N)$ (resp. $(\beta'_i~|~i\in\N)$)
the sequence of elementary divisors of $N$ (resp. of $N'$).
According to claim \ref{cl_directly-from-mouth} we have :
\set\begin{equation}\label{eq_rachmaninov}
\gamma_i-\beta_i,\gamma'_i-\beta'_i\in[-4\delta,4\delta]
\qquad
\text{for every $i\in\N$}.
\end{equation}
On the other hand, by the usual arguments we may find
$V$-linear maps $f:N\to M_*$ and $g:M'_*\to N'$ such that
$b^2\cdot\Coker\,f=b^2\cdot\Coker\,g=0$ for every $b\in V$
with $\log|b|_V=\delta$. Set $h:=g\circ\psi_*\circ f$;
it follows easily that $b^5\cdot\Coker\,h=0$. The sequence
of elementary divisors of $b^5N'$ is
$(\max(0,\beta'_i-5\delta)~|~i\in\N)$, and since
$N'':=\Img\,h$ is a finitely presented $V$-module, the
sequence $(\beta''_i~|~i\in\N)$ of its elementary divisors
satisfies the inequalities
$$
\beta_i\geq\beta''_i\geq\beta'_i-5\delta
\qquad
\text{for every $i\in\N$}
$$
by lemma \ref{lem_rachmaninov}(i,ii). Combining with
\eqref{eq_rachmaninov}, we deduce that
$$
\gamma_i\geq\gamma'_i-13\delta
\qquad
\text{for every $i\in\N$}
$$
and since $\delta$ is arbitrary, the assertion follows.

(ii): in light of \cite[Lemma 2.3.7(iii)]{Ga-Ra} and
claim \ref{cl_first-for-V-mods}(i), it is easily seen
that $\phi^\dagger$ is an epimorphism, so the assertion
follows formally from (i) and (iii), as in the proof of
lemma \ref{lem_rachmaninov}(i).
\end{proof}

One may ask, to which extent the sequence of elementary
divisors of a $V^a$-module $M$ determines the isomorphism
class of $M$. To address this question, we make the
following :

\begin{definition}\label{def_M-equivalence}
Let $(V,\fm)$ be an arbitrary basic setup, $A$ any
$V^a$-algebra, $M$ an $A$-module and $\omega$
a (very large) cardinal number, such that $\cM_\omega(A)$
contains the isomorphism class of $M$. The
topological space underlying the uniform space
$\cM_\omega(A)$ is not necessarily separated, but it admits
a maximal separated quotient space $\cM^\sep_\omega(A)$.
The {\em $\omega$-type of $M$} is the image in
$\cM^\sep_\omega(A)$ of the isomorphism class of $M$.
If $\omega'$ is any cardinal larger than $\omega$,
clearly the $\omega$-type of $M$ maps to the
$\omega'$-type of $M$, under the natural map
$\cM^\sep_\omega(A)\to\cM^\sep_{\omega'}(A)$. For
this reason, we shall usually omit explicit mention
of $\omega$, and shall call the {\em type of $M$}
any one of these invariants.
\end{definition}

\begin{lemma}\label{lem_Scholze}
In the situation of definition {\em\ref{def_M-equivalence}},
suppose that $M$ and $M'$ are two $A$-modules of the same
type. Then we have :
\begin{enumerate}
\item
$M$ is almost finitely presented (resp. almost finitely
generated) if and only if the same holds for $M'$.
\item
$M$ is a flat (resp. almost projective) $A$-module
if and only if the same holds for $M'$.
\end{enumerate}
\end{lemma}
\begin{proof} (i) follows directly from the definitions.

(ii) follows easily from \cite[Lemma 2.3.7(iii,iv)]{Ga-Ra} :
details left to the reader.
\end{proof}

We may now answer as follows to the foregoing question :

\begin{proposition}\label{prop_Scholze}
Let $(V,\fm_V)$ be as in \eqref{subsec_replace-cor}, and
$M$, $M'$ any two almost finitely generated torsion
$V^a$-modules. The following conditions are equivalent :
\begin{enumerate}
\alphaenu
\item
$M$ and $M'$ have the same type.
\item
The sequences of elementary divisors of $M$ and $M'$
coincide.
\end{enumerate}
\end{proposition}
\begin{proof}(a)$\Rightarrow$(b): Condition (a) means that
$(M,M')\in E_\delta$ for every
$\delta\in\log\Gamma^+_{\!\!V}\setminus\{0\}$, so the assertion is
immediate from the construction of the sequences of elementary
divisors.

(b)$\Rightarrow$(a): Let $(\gamma_i~|~i\in\N)$ be the common
sequence of elementary divisors for $M$ and $M'$. Fix $\delta$
as in the foregoing, pick finitely presented $V$-modules
$N$, $N'$ such that
$$
(M,N),(M',N')\in E_\delta
$$
and let $(\beta_i~|~i\in\N)$, respectively $(\beta'_i~|~i\in\N)$
be the sequences of elementary divisors of $N$ and $N'$.
From claim \ref{cl_directly-from-mouth} we deduce
$$
\beta_i-\gamma_,\beta'_i-\gamma_i\in[-4\delta,4\delta]
\qquad
\text{for every $i\in\N$}.
$$
Thus, $\beta_i-\beta'_i\in[-8\delta,8\delta]$ for every
$i\in\N$. Say that $8\delta=\log|a|_V$ for some $a\in V$.

\begin{claim}\label{cl_Scholze}
There exists a $V$-linear map $f:N\to N'$ whose kernel and
cokernel is annihilated by $a$.
\end{claim}
\begin{pfclaim} By \cite[Lemma 6.1.14]{Ga-Ra} we reduce
easily to the case where $N$ and $N'$ are cyclic $V$-modules,
so $N=V/bV$ and $N'=V/b'V$, for $b,b'\in V\setminus\{0\}$ such
that $\log|b/b'|_V\in[-2\delta,2\delta]$.  Now, if $b\in b'V$,
we let $f$ be the natural projection $V/bV\to V/b'V$.
If $b'\in bV$, we notice that $b'b^{-1}\one_N$ factors through
a map $N\to N'$ with the sought properties (details left to
the reader).
\end{pfclaim}

Clearly claim \ref{cl_Scholze} implies that $(N,N')\in E_{8\delta}$,
and therefore $(M,M')\in E_{10\delta}$. Since $\delta$ is
arbitrary, the assertion follows.
\end{proof}

\begin{corollary}\label{cor_flat-eldivs}
Let $b\in\fm_V$ be any non-zero element, $M$ an almost
finitely generated $V^a/bV^a$-module, and $(\gamma_i~|~i\in\N)$
the sequence of elementary divisors of $M$. The following
conditions are equivalent :
\begin{enumerate}
\alphaenu
\item
$M$ is a flat $V^a/bV^a$-module.
\item
There exists $n\in\N$ such that $\gamma_i=\log|b|_V$ for every
$i\leq n$, and $\gamma_i=0$ for every $i>n$.
\end{enumerate}
\end{corollary}
\begin{proof} For every $i\in\N$, pick $a_i\in V$ with
$\log|a_i|_V=\gamma_i$. According to proposition
\ref{prop_Scholze}, the type of $M$ coincides with that
of $M':=\bigoplus_{i\in\N}V^a/a_iV^a$, and by lemma
\ref{lem_Scholze}(ii), we may then assume that $M=M'$,
in which case clearly (b)$\Rightarrow$(a). For the
converse, notice that a standard computation gives
$$
\Tor_1^{V/bV}(V/aV,V/aV)\simeq V/(aV+a^{-1}bV)
\qquad
\text{for every $a\in V$ such that $\log|a|_V\leq\log|b|_V$}
$$
from which the assertion follows easily (details left
to the reader).
\end{proof}

\begin{proposition}
Let $(V,\fm_V)$ be as in \eqref{subsec_replace-cor}, and
$M$ a torsion $V$-module such that $M^a$ is almost finitely
generated. Let $(\gamma_i~|~i\in\N)$ be the sequence
of elementary divisors of $M^a$. Then
$$
\lambda_V(M)=(\Gamma_{\!V}:\Gamma)\cdot\sum_{i\in\N}\gamma_i.
$$
Especially, $\lambda_V(M)$ depends only on the type of $M^a$.
\end{proposition}
\begin{proof} Fix a sequence $(a_k~|~k\in\N)$ of elements of
$\fm_V\setminus\{0\}$ with $\lim_{k\to+\infty}\log|a_k|_V=0$.
For every $k\in\N$, let $M_k\subset M$ be a finitely generated
submodule such that $a_kM\subset M_k$
(\cite[Prop. 2.3.10(i)]{Ga-Ra}), and denote also by
$(\gamma_i^k~|~i\in\N)$ the sequence of elementary divisors
of $M_k^a$. After replacing each $M_k$ by $\sum_{j=0}^kM_j$, we
may also assume that $M_k\subset M_{k+1}$ for every $k\in\N$.

\begin{claim}\label{cl_gamma-bounds}
For every $i,k\in\N$ we have
$$
4\delta_k\geq\gamma_i-\gamma^k_i\geq 0
\qquad
\gamma^{k+1}_i\geq\gamma^k_i
\qquad
\text{where $\delta_k:=\log|a_k|$}.
$$
\end{claim}
\begin{pfclaim} The inequalities
$\gamma_i\geq\gamma_i^{k+1}\geq\gamma^k_i$
follow from proposition \ref{prop_rachmaninov}(ii).
Next, clearly we have $(M^a_k,M^a)\in E_{\delta_k}$ for
every $k\in\N$, so the upper bound for $\gamma_i-\gamma_i^k$
follows from claim \ref{cl_directly-from-mouth}.
\end{pfclaim}

Set $L^k:=\sum_{i\in\N}\gamma^k_i$ and
$L:=\sum_{i\in\N}\gamma_i$. It follows easily from claim
\ref{cl_gamma-bounds} that
$$
\lim_{k\to+\infty}L^k=\sum_{i\in\N}\lim_{k\to+\infty}\gamma^k_i=L.
$$
In light of proposition \ref{prop_additive-vr}(iii.b),
it then suffices to check the proposition with $M$
replaced by $M_k$, for every $k\in\N$. We may therefore
assume from start that $M$ is finitely generated, and
let $g$ be the cardinality of a finite system of generators
for $M$. We may then find a surjective $V$-linear map
$f:F:=V^{\oplus g}\to M$, and $N:=\Ker\,f^a$ is almost
finitely generated (claim \ref{cl_sub-of-uni}). Thus,
for every $k\in\N$ we may find a finitely generated
$V$-submodule $N_k\subset N$ with $a_kN\subset N_k$,
so that $(F^a/N_k^a,M)\in E_{\delta_k}$.
After replacing each $N_k$ by $\sum_{j=0}^kN_j$, we
may also assume that $N_k\subset N_{k+1}$ for every $k\in\N$.
Denote by $(\beta^k_i~|~i\in\N)$ the sequence of elementary
divisors of $F/N_k$, and notice that $\beta^k_i=0$ for every
$i\geq g$. Taking into account \eqref{eq_div-comp-length},
claim \ref{cl_directly-from-mouth} and proposition
\ref{prop_rachmaninov}(iii) we get
$$
\lambda_V(F/N_k)=\sum_{i\in\N}\beta^k_i
\qquad
4\delta_k\geq\beta_i^k-\gamma_i\geq 0
\qquad
\text{for every $i,k\in\N$}.
$$
Therefore $\lim_{k\to+\infty}\lambda_V(F/N_k)=L$,
and in view of proposition \ref{prop_additive-vr}(i),
it remains only to check that
$\lim_{k\to+\infty}\lambda_V(N/N_k)=0$, or equivalently,
that $\lim_{k\to+\infty}\lambda_V(N_k)=\lambda_V(N)$
(proposition \ref{prop_additive-vr}(i)). However,
set $N':=\bigcup_{k\in\N}N_k$; by proposition
\ref{prop_additive-vr}(iii.b), the latter identity
holds if and only $\lambda_V(N')=\lambda_V(N)$,
and since $(N/N')^a=0$, this follows from proposition
\ref{prop_additive-vr}(i,ii).
\end{proof}

\sset\subsubsection{}\label{subsec_one-more-cat}
In order to deal with general measurable $K^+$-algebras, we
introduce hereafter some further notation which shall be
standing throughout this section.

$\bullet$\ \
To begin with, any ring homomorphism $\phi:A\to B$ induces
functors
\set\begin{equation}\label{eq_restr-scalars}
\phi_*:B\Mod\to A\Mod
\qquad\text{and}\qquad
\phi^*:A\Mod\to B\Mod.
\end{equation}
Namely, $\phi_*$ assigns to any $B$-module $M$ the $A$-module
$\phi_*M$ obtained by restriction of scalars, and
$\phi^*(M):=B\otimes_AM$.

$\bullet$\ \
For any local ring $(A,\fm_A)$, we let $\kappa(A):=A/\fm_A$,
and we denote by $s(A)$ the closed point of $\Spec\,A$.
If $A$ is also a coherent ring, we denote by $A\Mod_\cohs$
the full subcategory of $A\Mod$ consisting of all the finitely
presented $A$-modules $M$ such that $\Supp\,M\subset\{s(A)\}$.
Notice that the coherence of $A$ implies that $A\Mod_\cohs$ is
an abelian category.

$\bullet$\ \
Lastly, let $\cA$ be any small abelian category; recall
that $K_0(\cA)$ is the abelian group defined by generators
and relations as follows. The generators are the isomorphism
classes $[T]$ of objects $T$ of $\cA$, and the relations are
generated by the elements of the form $[T_1]-[T_2]+[T_3]$,
for every short exact sequence $0\to T_1\to T_2\to T_3\to 0$
of objects of $\cA$. One denotes by
$K_0^+(\cA)\subset K_0(\cA)$ the submonoid generated by the
classes $[T]$ of all objects of $\cA$.
We shall use the following well known {\em d\'evissage\/}
lemma :

\begin{lemma}\label{lem_devissage}
Let $\iota:\cB\subset\cA$ be an additive exact and fully
faithful inclusion of abelian categories, and suppose that :
\begin{enumerate}
\alphaenu
\item
If\/ $T\in\Ob(\cB)$ and $T'$ is a subquotient of $\iota(T)$,
then $T'$ is in the essential image of $\iota$.
\item
Every object $T$ of $\cA$ admits a finite filtration
$\Fil^\bullet T$ such that the associated graded object
$\gr^\bullet T$ is in the essential image of $\iota$.
\romanenu
\end{enumerate}
Then $\iota$ induces an isomorphism :
$$
K_0(\cB)\isom K_0(\cA).
$$
\end{lemma}
\begin{proof} Left to the reader.
\end{proof}

\begin{proposition}\label{prop_about_K_0}
Let $\phi:A\to B$ be a morphism of measurable $K^+$-algebras. We have :
\begin{enumerate}
\item
If $\phi$ induces a finite field extension $\kappa(A)\to\kappa(B)$,
then the functor $\phi_*$ of \eqref{eq_restr-scalars} restricts to
a functor
$$
\phi_*:B\Mod_\cohs\to A\Mod_\cohs
$$
which induces a group homomorphism of the respective
$K_0$-groups :
$$
\phi_*:K_0(B\Mod_\cohs)\to K_0(A\Mod_\cohs).
$$
\item
If $\phi$ induces an integral morphism $\kappa(A)\to B/\fm_AB$,
then $\length_B(B/\fm_AB)$ is finite, and the functor $\phi^*$
of \eqref{eq_restr-scalars} restricts to a functor
$$
\phi^*:A\Mod_\cohs\to B\Mod_\cohs
$$
and if $\phi$ is also a flat morphism, $\phi^*$ induces a group
homomorphism :
$$
\phi^*:K_0(A\Mod_\cohs)\to K_0(B\Mod_\cohs).
$$\end{enumerate}
\end{proposition}
\begin{proof} Write $A$ (resp. $B$) as the colimit of
a filtered system $\underline A:=(A_i~|~i\in I)$ (resp.
$\underline B:=(B_j~|~j\in J)$) of objects of $K^+\mAlg_0$,
with local and essentially \'etale transition maps. After
replacing $J$ (resp. $I$) by a cofinal subsets, we may
assume that the indexing set admits an initial element
$0\in J$ (resp. $0\in I$). Furthermore, we may assume
that the induced map $A_0\to B$ factors through a morphism
$A_0\to B_0$ in $K^+\mAlg_0$.

(i): Let $M$ be any object of $B\Mod_\cohs$. We
need to show that $\phi_*M$ is finitely presented. 
We may find $j\in J$ and a finitely presented $B_j$-module
$M_j$, with an isomorphism $M\isom M_j\otimes_{B_j}B$
of $B$-modules. After replacing $J$ by $J/j$, we may assume
that $j=0$ is the initial index. Since the natural map
$B_0\to B$ is local and ind-\'etale, it is easily seen that
$M_0$ is an object of $B_0\Mod_\cohs$. Especially, there
exists a finitely generated $\fm_{B_0}$-primary ideal
$I\subset\Ann_{B_0}M_0$. We may then replace the system
$\underline B$ by $(B_j/IB_j~|~j\in J)$ and assume that
each $B_j$ has Krull dimension zero.

\begin{claim}\label{cl_find}
Let $\phi:A\to B$ be a morphism of measurable $K^+$-algebras
inducing a finite residue field extension $\kappa(A)\to\kappa(B)$,
and such that $B$ has dimension zero. Let also $M$
be any finitely presented $B$-module. Then we may
find :
\begin{enumerate}
\alphaenu
\item
a cocartesian diagram of local maps of $K^+$-algebras
$$
\xymatrix{ A_l \ar[r]^-{\phi_l} \ar[d] & C_l \ar[d] \\
           A \ar[r]^-\phi & B
}$$
whose vertical arrows are ind-\'etale, and where $\phi_l$
is a morphism in $K^+\mAlg_0$
\item
and an object $M_l$ of $C_l\Mod_\cohs$, with an isomorphism
$M_l\otimes_{C_l}B\isom M$.
\end{enumerate}
\end{claim}
\begin{pfclaim} Define $\underline A$ and $\underline B$ as in
the foregoing. Notice that -- under the current assumptions --
$B_j$ is a henselian ring, for every $j\in J$. Moreover, we may
find $j\in J$ such that $\kappa(B)$ is generated by the image
of $\kappa(A)\otimes_{\kappa(B_0)}\kappa(B_j)$; after
replacing again $J$ by $J/j$, we may then also assume
that $\kappa(B)=\kappa(A)\cdot\kappa(B_0)$.

For every $i\in I$, set $B'_i:=A_i\otimes_{A_0}B_0$; the
natural map $B_0\to B$ factors through a map $B'_i\to B$,
and we let $\fp_i\subset B'_i$ be the preimage of $\fm_B$.
Set also $C_i:=B'_{i,\fp_i}$, so we deduce a filtered
system of local maps $(C_i\to B~|~i\in I)$, whose limit
is a local map $\psi:C\to B$ of local ind-\'etale
$B_0$-algebras, which -- by construction -- induces an
isomorphism $\kappa(C)\isom\kappa(B)$ on residue fields.
It follows easily that $\psi$ is itself ind-\'etale, so
say that $\psi$ is the colimit of a filtered system
$(\psi_\lambda:C\to D_\lambda~|~\lambda\in\Lambda)$ of
\'etale $C$-algebras.
Notice that $C$ is a henselian local ring; in light of
\cite[Ch.IV, Th.18.5.11]{EGA4} we may then assume that
$D_\lambda$ is a local ring and $\psi_\lambda$ is a finite
\'etale map, for every $\lambda\in\Lambda$. Clearly the induced
residue field extension $\kappa(C)\to\kappa(D_\lambda)$
is an isomorphism; in view of \cite[Ch.IV, Prop.18.5.15]{EGA4}
it follows that $\psi_\lambda$ is an isomorphism, for every
$\lambda\in\Lambda$, so the same holds for $\psi$.

Notice that the sequence of residue degrees
$d_i:=[\kappa(C_i):\kappa(A_i)]$ is non-increasing,
hence there exists $l\in I$ such that $d_i=d:=d_l$
for every index $i\geq l$. Notice as well that, for
$i\geq l$, the local algebra $C_i$ is also a localization
of $C'_i:=A_i\otimes_{A_l}C_l$, and the latter is an
essentially finitely presented $K^+$-algebras of Krull
dimension zero, hence its spectrum is finite and discrete
(lemma \ref{lem_fin-min}). Moreover, since  the image of
the map $\Spec\,C_l\to\Spec\,A_l$ is the closed point,
it is clear that the same holds for the image of the
induced map $\Spec\,C'_i\to\Spec\,A_i$. Since the
extension $\kappa(A_l)\to\kappa(A_i)$ is finite and
separable, we conclude that
\set\begin{equation}\label{eq_these-are-them}
\kappa(A_i)\otimes_{\kappa(A_l)}\kappa(C_l)=\!\!\!\!
\prod_{\fp\in\Spec\,C'_i}\!\!\!\kappa(C'_{i,\fp}).
\end{equation}
However, clearly the left-hand side of \eqref{eq_these-are-them}
is a $\kappa(A_i)$-algebra of degree $d$, whereas one
of factors of the right-hand side -- namely $\kappa(C_i)$ --
is already of degree $d$ over $\kappa(A_i)$. Hence
$\Spec\,C'_i$ contains a single element, {\em i.e.}
$C'_i=C_i$ is a local ring, and $C=A\otimes_{A_l}C_l$.
Summing up, we have obtained the sought cocartesian
diagram, and the claim holds with $M_l:=M_0\otimes_{B_0}C_l$.
\end{pfclaim}

Let $M_l$ and $\phi_l$ be as in claim \ref{cl_find}; then
$\phi_*M$ is isomorphic to $A\otimes_{A_l}\phi_{l*}M_l$, so
may replace from start $\phi$ by $\phi_l$, and assume that
$\phi$ is a morphism in $K^+\mAlg_0$, with $B$ of Krull
dimension zero. In such situation, one sees easily that
$\phi$ is integral, hence $B/I$ is a finitely presented
$A$-module, by proposition \ref{prop_integral-fp-ext}(i),
therefore $\phi_*M$ is a finitely presented $A$-module, as
required. Lastly, since the functor $\phi_*$ is exact, it
is clear that it induces a map on $K_0$-groups as stated.

(ii): Under the current assumptions, the induced map
$\kappa(A_0)\to B_0/\fm_{A_0}B_0$ is integral and essentially
finitely presented, hence it is finite, so $B_0/\fm_{A_0}B_0$
is a $B_0$-module of finite length; but this is also the
length of the $B$-module $B/\fm_AB$, whence the first
assertion. Next, let $M$ be an object of $A\Mod_\cohs$; then
$\phi^*M$ is a finitely presented $B$-module; moreover, we
may find a $\fm_A$-primary ideal $I\subset A$ such that $M$
is a $A/I$-module, hence $\phi^*M$ is a $B/IB$-module. Notice
that the induced map $\Spec\,B/\fm_AB\to\Spec\,B/IB$ is bijective,
and its target is a local scheme of dimension zero (since
$B/\fm_AB$ is integral over a field). It follows easily that
$IB$ is a $\fm_B$-primary ideal, so $\phi^*M$ is an object of
$B\Mod_\cohs$. The last assertion is then a trivial consequence
of the exactness of the functor $\phi^*$, when $\phi$ is flat.
\end{proof}

\begin{lemma}\label{lem_find-val-ring}
Let $A$ be any measurable $K^+$-algebra. Then there exists a
morphism
$$
V \xrightarrow{\ \phi\ }A/I
$$
of measurable $K^+$-algebras, where :
\begin{enumerate}
\alphaenu
\item
$I\subset A$ is a finitely generated $\fm_A$-primary ideal.
\item
$V$ is a valuation ring, $\phi$ is a finitely presented
surjection and the natural map $K^+\to V$ induces an
isomorphism of value groups $\Gamma\isom\Gamma_{\!\!V}$.
\end{enumerate}
\end{lemma}
\begin{proof} Suppose first that $A$ is an object of
$K^+\mAlg_0$. In this case, choose an affine finitely
presented $S$-scheme $X$ and a point $x\in X$ such that
$A=\cO_{\!X,x}$; next, take a finitely presented closed
immersion $h:X\to Y:=\A^n_{K^+}$ of $S$-schemes; set
$\bar Y:=\A^n_\kappa\subset Y$, pick elements
$f_1,\dots,f_d\in B:=\cO_{Y,h(x)}$ whose images in the
regular local ring $\cO_{\bar Y,h(x)}$ form a regular
system of parameters ({\em i.e.} a regular sequence that
generates the maximal ideal), and let $J\subset B$ be the
ideal generated by the $f_i$, $i=1,\dots,d$. Let $I\subset A$
be any finitely generated $\fm_A$-primary ideal containing the
image of $J$. We deduce a surjection $\phi:V:=B/J\to A/I$, and
by construction $V/\fm_K V$ is a field; moreover, the induced
map $K^+\to V$ is flat by virtue of \cite[Ch.IV, Th.11.3.8]{EGAIV-3}.
It follows that $V$ is a valuation ring with the sought
properties, by proposition \ref{prop_integral-fp-ext}(ii).
Also, $\phi$ is finitely presented, by proposition
\ref{prop_integral-fp-ext}(i).

Next, let $A$ be a general measurable $K^+$-algebra, and write
$A$ as the colimit of a filtered system $(A_j~|~j\in J)$ of
objects of $K^+\mAlg_0$. We may assume that $0\in J$ is an
initial index, and the foregoing case yields an
$\fm_{A_0}$-primary ideal $I_0\subset A_0$, and a surjective
finitely presented morphism $\phi_0:V_0\to A_0/I_0$ in
$K^+\mAlg_0$ from a valuation ring $V_0$, such that
$\Gamma_{\!\!V}=\Gamma$. Notice that $A/I_0$ is a henselian
ring, hence $\phi_0$ extends to a ring homomorphism
$\phi_0^\he:V_0^\he\to A_0/I_0$ from the henselization
$V_0^\he$ of $V_0$; more precisely, $\phi_0$ induces an
isomorphism of $V^\he_0$-algebras :
$$
V_0^\he\otimes_{V_0}(A_0/I_0)\isom A_0/I_0
$$
so $\phi_0^\he$ is still finitely presented. On the one hand,
$\phi_0^\he$ induces an identification
$$
\kappa(V^\he_0)=\kappa(A_0).
$$
On the other hand, we have the filtered system of separable
field extensions $(\kappa(A_j)~|~j\in J)$, whose colimit is
$\kappa(A)$. There follows a corresponding filtered system
$(V^\he_j~|~j\in J)$ of finite \'etale $V_0^\he$-algebras,
whose colimit we denote $V$ (\cite[Ch.IV, Prop.18.5.15]{EGA4}).
Then $V$ is a valuation ring, and the map $V^\he_0\to V$
induces an isomorphism on value groups. Moreover, the
induced isomorphisms $\kappa(V^\he_j)\isom\kappa(A_j)$
lift uniquely to morphisms of $A_0$-algebras
$\phi^\he_j:V^\he_j\to A_j/I_0A_j$, for every $j\in J$
(\cite[Ch.IV, Cor.18.5.12]{EGA4}). Due to the uniqueness
of $\phi_j^\he$, we see that the resulting system
$(\phi_j^\he~|~j\in J)$ is filtered, and its colimit
is a morphism $\phi:V\to A/I_0A$. Moreover, $\phi_j^\he$
induces an isomorphism
$V^\he_j\otimes_{V_0^\he}A_0/IA_0\isom A_j/I_0A_j$, especially
$\phi_j^\he$ is surjective for every $j\in J$, so the same
holds for $\phi$. More precisely, $\phi_0^\he$ induces an
isomorphism $V\otimes_{V^\he_0}(A_0/I_0)\isom A/I_0A$, hence
$\phi$ is still finitely presented.
\end{proof}

\begin{theorem}\label{th_about_K_0}
With the notation of \eqref{subsec_one-more-cat}, the following
holds :
\begin{enumerate}
\item
For every measurable $K^+$-algebra $A$ there is a natural
group isomorphism :
$$
\blambda_A:K_0(A\Mod_\cohs)\isom\log\Gamma
$$
which induces an isomorphism $K^+_0(A\Mod_\cohs)\isom\log\Gamma^+$.
\item
The family of isomorphisms $\blambda_A$ (for $A$ ranging over
the measurable $K^+$-algebras) is characterized uniquely
by the following two properties.
\begin{enumerate}
\item
If\/ $V$ is a valuation ring and a flat measurable $K^+$-algebra,
then
$$
\blambda_V([M])=\lambda_V(M)
\qquad
\text{for every object $M$ of\/ $V\Mod_\cohs$}
$$
where $\lambda_V(M)$ is defined as in \eqref{subsec_aim-length}.
\item
Let\/ $\psi:A\to B$ be a morphism of measurable $K^+$-algebras
inducing a finite residue field extension $\kappa(A)\to\kappa(B)$.
Then
$$
\blambda_A(\psi_*[M])=[\kappa(B):\kappa(A)]\cdot\blambda_B([M])
\qquad
\text{for every $[M]\in K_0(B\Mod_\cohs)$}.
$$
\end{enumerate}
\item
For every $a\in K^+\!\setminus\!\{0\}$ and any object $M$ of
$A\Mod_\cohs$ we have :
\begin{enumerate}
\item
$[M]=0$ in $K_0(A\Mod_\cohs)$ if and only if $M=0$.
\item
If $M$ is flat over $K^+/aK^+$, then :
\set\begin{equation}\label{eq_charct-blambda}
\blambda_A([M])=|a|\cdot\length_A(M\otimes_{K^+}\kappa).
\end{equation}
\end{enumerate}
\item
Let $\psi:A\to B$ be a flat morphism of measurable $K^+$-algebras
inducing an integral map $\kappa(A)\to B/\fm_AB$. Then :
$$
\blambda_B(\psi^*[M])=\length_B(B/\fm_A B)\cdot\blambda_A([M])
\qquad
\text{for every $[M]\in K_0(A\Mod_\cohs)$}.
$$
\end{enumerate}
\end{theorem}
\begin{proof} We start out with the following :

\begin{claim}\label{cl_bottoms-up}
Let $(K,|\cdot|)\to(E,|\cdot|)$ be an extension of valued
fields of rank one inducing an isomorphism of value groups,
$a\in K^+\setminus\{0\}$ any element, $M$ an $E^+/aE^+$ module.
Then :
\begin{enumerate}
\item
$M$ is a flat $E^+/aE^+$-module if and only if it is a
flat $K^+/aK^+$-module.
\item
$M\otimes_{K^+}\kappa=M\otimes_{E^+}\kappa(E)$.
\end{enumerate}
\end{claim}
\begin{pfclaim} According to \cite[Th.7.8]{Mat}, in order to
show (i) it suffices to prove that
$$
\Tor^{E^+/aE^+}_1(E^+/bE^+,M)=\Tor^{K^+/aK^+}_1(K^+/bK^+,M)
$$
for every $b\in K^+$ such that $|b|\geq|a|$. The latter
assertion is an easy consequence of the faithful flatness
of the extension $K^+\to E^+$.
(ii) follows from the identity : $\fm_E=\fm_K E^+$, which holds
since $\Gamma_{\!\!E}=\Gamma$.
\end{pfclaim}

\begin{claim}\label{cl_OK-for-vals}
Let $f:K^+\to V$ be a morphism of measurable $K^+$-algebras,
where $V$ is a valuation ring and $f$ induces an isomorphism
on value groups. Then $\blambda_V$ (defined by (ii.a)) is an
isomorphism and assertion (iii) holds for $A=V$.
\end{claim}
\begin{pfclaim} According to (ii.a), $\blambda_V([M])=\lambda_V(M)$
for every finitely presented torsion $V$-module. However,
every such module $M$ admits a decomposition of the form
$M\simeq(V/a_1V)\oplus\cdots\oplus(V/a_kV)$, with
$a_1,\dots,a_k\in\fm_K$ (\cite[Lemma 6.1.14]{Ga-Ra}). Then
by claim \ref{cl_bottoms-up}(i), $M$ is flat over $K^+/aK^+$
if and only if $M$ is flat over $V/aV$, if and only if $|a|=|a_i|$
for every $i\leq k$. In this case, an explicit calculation
shows that $\lambda_V(M)=|a|\cdot\length(M\otimes_V\kappa(V))$,
which is equivalent to \eqref{eq_charct-blambda}, in view of
claim \ref{cl_bottoms-up}(ii). Next, we consider the map :
$$
\bmu:\log\Gamma^+\to K_0(V\Mod_\cohs)
\quad :\quad
|a|\mapsto[V/aV]
\quad\text{for every $a\in K^+\setminus\{0\}$}.
$$
We leave to the reader the verification that $\bmu$
extends to a group homomorphism well-defined on the whole
of $\Gamma$, that provides an inverse to $\blambda_V$.
Finally, it is clear that $\lambda_V(M)=0$ if and only if
$M=0$, so also (iii.a) holds.
\end{pfclaim}

\begin{claim}\label{cl_filter-that}
Let $A$ be any measurable $K^+$-algebra. For every object
$N$ of $A\Mod_\cohs$ there exist a finite filtration
$0=N_0\subset\cdots\subset N_k=N$ by finitely presented
$A$-submodules, and elements $a_1,\dots,a_k\in\fm$ such that
$N_i/N_{i-1}$ is a flat $K^+/a_iK^+$-module for every
$1\leq i\leq k$.
\end{claim}
\begin{pfclaim} Let us find $I\subset A$ and $\phi:V\to A/I$ as
in lemma \ref{lem_find-val-ring}. It suffices to show the
claim for the finitely presented $A$-modules
$I^nN/I^{n+1}N$ (for every $n\in\N$), hence we may assume that
$N$ is an $A/I$-module. Then $\phi_*N$ is a finitely presented
$V$-module, hence of the form $(V/a_1V)\oplus\cdots\oplus(V/a_kV)$
for some $a_i\in\fm$; we may order the summands so that
$|a_i|\geq|a_{i+1}|$ for all $i<k$. We argue by induction on
$d(N):=\dim_\kappa(N\otimes_V\kappa(V))$. If $d(N)=0$, then
$N=0$ by Nakayama's lemma. Suppose $d>0$; we remark that
$N/a_1N$ is a flat $V/a_1V$-module, hence a flat
$K^+/a_1K^+$-module (claim \ref{cl_bottoms-up}),
and $d(a_1N)<d(N)$; the claim follows.
\end{pfclaim}

\begin{claim}\label{cl_devis-A-Mod}
Let $A$ be any measurable $K^+$-algebra, $I\subset A$ a finitely
generated $\fm_A$-primary ideal, and $\pi_I:A\to A/I$ the
natural projection. Then the map :
$$
\pi_{I*}:K_0(A/I\Mod_\coh)\to K_0(A\Mod_\cohs)
$$
is an isomorphism.
\end{claim}
\begin{pfclaim} On the one hand we have :
$$
A\Mod_\cohs=\bigcup_{n\in\N} A/I^n\Mod_\coh
$$
and on the other hand, in view of lemma \ref{lem_devissage},
we see that the projections $A/I^{n+1}\to A/I^n$ induce
isomorphisms $K_0(A/I^n\Mod_\coh)\to K_0(A/I^{n+1}\Mod_\coh)$
for every $n>0$, whence the claim.
\end{pfclaim}

Let $A$, $I$ and $\phi:V\to A/I$ be as in lemma
\ref{lem_find-val-ring}, and $\pi_I:A\to A/I$ the natural
surjection; taking into account claim \ref{cl_devis-A-Mod},
we may let :
\set\begin{equation}\label{eq_rule-for-blambda}
\blambda_A:=\blambda_V\circ\phi_*\circ\pi_{I*}^{-1}
\end{equation}
where $\blambda_V$ is given by the rule of (ii.a). In view of
claim \ref{cl_OK-for-vals} we see that $\blambda_A([M])=0$
if and only if $M=0$, so (iii.a) follows already.

\begin{claim}\label{cl_independent}
The isomorphism $\blambda_A$ is independent of the choice
of $I$, $V$ and $\phi$.
\end{claim}
\begin{pfclaim} Indeed, suppose that
$J\subset A$ is another ideal and $\psi:W\to A/J$ is another
surjection from a valuation ring $W$, fulfilling the
foregoing conditions. We consider the commutative diagram
$$
\xymatrix{ & A/I \ar[dr]^{\bar\pi_J} & V \ar[l]_-\phi \ar[d]^{\phi'} \\
A \ar[ur]^{\pi_I} \ar[dr]_{\pi_J} \ar[rr]^-{\pi_{I+J}} & & A/(I+J) \\
& A/J \ar[ru]_{\bar\pi_I} & W \ar[l]^\psi \ar[u]_{\psi'}.
}$$
We compute : $\phi_*\circ\pi_{I*}^{-1}=
\phi_*\circ\pi_{J*}\circ\bar\pi{}_{J*}^{-1}\circ\bar\pi{}_{I*}^{-1}=
\phi'_*\circ\pi_{I+J,*}^{-1}$, and a similar calculation
shows that
$\psi_*\circ\pi_{J*}^{-1}=\psi'_*\circ\pi_{I+J,*}^{-1}$.
We are thus reduced to showing that
$\blambda_V\circ\phi'_*([N])=
\blambda_W\circ\psi'_*([N])$ for every $A/(I+J)$-module $N$.
In view of claim \ref{cl_filter-that} we may assume that $N$ is flat
over $K^+/aK^+$, for some $a\in\fm$, in which case the assertion
follows from claim \ref{cl_OK-for-vals}.
\end{pfclaim}

Claims \ref{cl_OK-for-vals} and \ref{cl_independent} show
already that (i) holds. Next, let $\psi:A\to B$ be as in (ii.b).
Choose an ideal $I\subset A$ and a surjection $\phi:V\to A/I$ as
in lemma \ref{lem_find-val-ring}.
Let also $J\subset B$ be a finitely generated  $\fm_B$-primary
ideal containing $\psi(I)$, and $\bar\psi:A/I\to B/J$ the
induced map.

By inspecting the definitions we see that (ii.b) amounts
to the identity : $\blambda_V((\bar\psi\circ\phi)_*[M])=
[\kappa(B):\kappa(A)]\cdot\blambda_{B/J}([M])$ for
every finitely presented $B/J$-module $M$.
Furthermore, up to enlarging $J$, we may assume that there
is a finitely presented surjection $\xi:W\to B/J$ of
$K^+$-algebras, where $W$ is a valuation ring with value
group $\Gamma$, and then we come down to showing:
$$
\blambda_V((\bar\psi\circ\phi)_*[M])=
[\kappa(B):\kappa(A)]\cdot\blambda_W(\xi_*[M]).
$$
In view of claim \ref{cl_filter-that} we may also assume that
$M$ is a flat $K^+/aK^+$-module for some $a\in\fm$, in which
case the identity becomes :
$$
|a|\cdot\length_V(M\otimes_{K^+}\kappa)=[\kappa(B):\kappa(A)]
\cdot|a|\cdot\length_W(M\otimes_{K^+}\kappa)
$$
thanks to claim \ref{cl_OK-for-vals}. However, the latter is
an easy consequence of the identities : $\kappa(A)=\kappa(V)$
and $\kappa(B)=\kappa(W)$. Next, we show that (iii.b) holds
for a general measurable $K^+$-algebra $A$. Indeed, let $M$
be any object of $A\Mod_\cohs$; as usual, we may find a local
ind-\'etale map $A_0\to A$ from some object $A_0$ of
$K^+\mAlg_0$, and an object $M_0$ of $A_0\Mod_\cohs$ with
an isomorphism of $A$-modules $A\otimes_{A_0}M_0\isom M$
(cp. the proof of proposition \ref{prop_about_K_0}(i)).
Since (ii.b) is already proved, we have
$$
\blambda_A([M])=\blambda_{A_0}([M_0])
\quad\text{and}\quad
\length_A(M\otimes_{K^+}\kappa)=
\length_{A_0}(M\otimes_{K^+}\kappa).
$$
Therefore, we may replace $A$ by $A_0$, and assume that
$A$ is an object of $K^+\mAlg_0$. In this case,
by lemma \ref{lem_integral-fp-ext} we may find morphisms
$f:K^+\to V$ and $g:V\to A$ in $K^+\mAlg_0$
such that $V$ is a valuation ring with value group $\Gamma$
and $g$ induces a finite extension of residue fields
$\kappa(V)\to\kappa(A)$.
Let $N$ be an $A$-module supported at $s(A)$ and flat over
$K^+/aK^+$; we may find a finitely generated $\fm_A$-primary
ideal $I\subset A$ such that $I\subset\Ann_A(N)$, and since
(ii.b) is already known in general, we reduce to showing that
(iii.b) holds for the $A/I$-module $N$. However, the induced
map $\bar g:V\to A/I$ is finite and finitely presented
(proposition \ref{prop_integral-fp-ext}(i)), so another
application of (ii.b) reduces to showing that
\eqref{eq_charct-blambda} holds for $A:=V$ and $M:=\bar g_*N$,
in which case the assertion is already known by claim
\ref{cl_OK-for-vals}.

(ii.a): Suppose that $A$ is a valuation ring, and let $M$
be any object of $A\Mod_\cohs$. The sought assertion is obvious
when $\Gamma_{\!\!A}=\Gamma$, since in that case we can choose
$A=V$ and $\phi=\pi_I$ in \eqref{eq_rule-for-blambda}.
However, we know already that the rank of $A$ equals one, and
$e:=(\Gamma_{\!\!A}:\Gamma)$ is finite (see
\eqref{subsec_replace-cor}); we can then assume that $e>1$,
in which case corollary \ref{cor_integral-fp-ext}(ii) implies
that $\Gamma\simeq\Z$. Then it suffices to check (ii.a) for
$M=\kappa(A)$, which is a (flat) $\kappa$-module, so that --
by assertion (iii) -- one has $\blambda_A([M])=e\cdot\gamma_0$,
where $\gamma_0\in\log\Gamma_{\!\!A}^+$ is the positive generator.
In view of remark \ref{rem_shallbeleft}, we see that this value
agrees with $\lambda_V(M)$, as stated.

(iv): In view of claim \ref{cl_filter-that}, we may assume
that $M$ is a flat $K^+/aK^+$-module, for some $a\in\fm$, and
then the same holds for $\psi^*M$, since $\psi$ is flat.
In view of (iii.b), it then suffices to show that :
$$
\length_B(B\otimes_AM)=\length_A(M)\cdot\length_B(B/\fm_A B)
$$
for any $A$-module $M$ of finite length. In turn, this is easily
reduced to the case where $M=\kappa(A)$, for which the
identity is obvious.
\end{proof}

\sset\subsubsection{}\label{subsec_C_M}
Let $A$ be a measurable $K^+$-algebra. The next step consists
in extending the definition of $\blambda_A$ to the category
$A\Mod_\bsb$ of arbitrary $A$-modules $M$ supported at $s(A)$.
First of all, suppose that $M$ is finitely generated. Let
$\cC_M$ be the set of isomorphism classes of objects $M'$ of
$A\Mod_\cohs$ that admit a surjection $M'\to M$. Then we set :
$$
\lambda^*_A(M):=\inf\,\{\blambda_A([M'])~|~M'\in\cC_M\}
\in\log\Gamma^\wedge.
$$
Notice that -- by the positivity property of theorem
\ref{th_about_K_0}(i) -- we have $\lambda^*_A(M)=\blambda_A([M])$
whenever $M$ is finitely presented. Next, for a general object
of $A\Mod_\bsb$ we let :
\set\begin{equation}\label{eq_really-full}
\lambda_A(M):=\sup\,\{\lambda^*_A(M')~|~M'\subset M
\text{ and $M'$ is finitely generated}\}\in
\log\Gamma^\wedge\cup\{+\infty\}.
\end{equation}

\begin{lemma} If $M$ is finitely generated, then
$\lambda^*_A(M)=\lambda_A(M)$.
\end{lemma}
\begin{proof} Let $M'\subset M$ be a finitely generated
submodule; we have to show that $\lambda^*(M')\leq\lambda^*(M)$.
To this aim, let $f:N\to M$ and $g:N'\to M'$ be two surjections
of $A$-modules with $N\in\cC_M$ and $N'\in\cC_{M'}$; by filtering
the kernel of $f$ by the system of its finitely generated
submodules, we obtain a filtered system $(N_i~|~i\in I)$
of finitely presented quotients of $N$, with surjective
transition maps, such that $\colim_{i\in I}N_i=M$.  By
\cite[Prop.2.3.16(ii)]{Ga-Ra} the induced map $N'\to M$
lifts to a map $h:N'\to N_i$ for some $i\in I$. Since $A$
is coherent, $h(N')$ is a finitely presented $A$-module
with a surjection $h(N')\to M'$, hence
$\lambda^*_A(M')\leq\blambda_A([h(N')])\leq\blambda_A([N_i])
\leq\blambda_A(N)$. Since $N$ is arbitrary, the claim
follows.
\end{proof}

\begin{proposition}\label{prop_slight-gen}
{\em (i)}\ \ If $A$ is a valuation ring, \eqref{eq_really-full}
agrees with \eqref{eq_full-gener}.
\begin{enumerate}
\addenu
\item
If $(M_i~|~i\in I)$ is a filtered system of objects of
$A\Mod_\bsb$ with injective (resp. surjective) transition
maps, then :
$$
\lambda_A(\colim_{i\in I}M_i)=\lim_{i\in I}\,\lambda_A(M_i)
$$
(resp. provided there exists $i\in I$ such that
$\lambda_A(M_i)<+\infty$).
\item
If\/ $\psi:A\to B$ is a morphism of measurable $K^+$-algebras
inducing a finite extension $\kappa(A)\to\kappa(B)$ of residue
fields, then :
$$
\lambda_B(M)=
\frac{\lambda_A(\psi_*M)}{[\kappa(B):\kappa(A)]}
\qquad
\text{for every object $M$ of $B\Mod_\bsb$}.
$$
\item
Let\/ $0\to M_1\to M_2\to M_3\to 0$ be a short exact
sequence in $A\Mod_\bsb$. Then :
$$
\lambda_A(M_2)=\lambda_A(M_1)+\lambda_A(M_3).
$$
\item
Let $\psi$ be a flat morphism of measurable $K^+$-algebras
inducing an integral map $\kappa(A)\to B/\fm_AB$. Then :
$$
\lambda_B(\psi^*M)=\length_B(B/\fm_A B)\cdot\lambda_A(M)
\qquad
\text{for every object $M$ of $A\Mod_\bsb$}.
$$
\end{enumerate}
\end{proposition}
\begin{proof} (i): Set $e:=(\Gamma_{\!\!V}:\Gamma)$.
It suffices to check the assertion for a finitely generated
$A$-module $M$, in which case one has to show the identity :
$$
e\cdot|F_0(M)^a|=\lambda_A^*(M).
$$
By proposition \ref{prop_additive-vr}(i) we have
$|F_0(M')^a|\geq|F_0(M)^a|$ for every $M'\in\cC_M$, hence
$\lambda^*_A(M)\geq e\cdot|F_0(M)^a|$.
On the other hand, let us fix a surjection $V^{\oplus n}\to M$,
and let us write its kernel in the form
$K=\bigcup_{i\in I}K_i$, for a filtered system $(K_i~|~i\in I)$
of a finitely generated $V$-submodules; it follows that :
$$
F_0(M)=\bigcup_{i\in I}F_0(V^{\oplus n}/K_i).
$$
In view of proposition \ref{prop_additive-vr}(iii.b)
we deduce :
$e\cdot|F_0(M)^a|=\lim_{i\in I}\,e\cdot|F_0(V^{\oplus n}/K_i)^a|
\geq\lambda_A^*(M)$.

The proof of (ii) in the case where the transition
maps are injective, is the same as that of proposition
\ref{prop_additive-vr}(iii.b).

(iii): Let us write $M=\bigcup_{i\in I}M_i$ for a filtered
family $(M_i~|~i\in I)$ of finitely generated $B$-submodules.
By the case already known of (ii) we have :
$\lambda_B(M)=\lim_{i\in I}\,\lambda_B(M_i)$, and likewise
for $\lambda_A(\psi_*M)$, hence we may assume
from start that $M$ is a finitely generated $B$-module, in
which case the annihilator of $M$ contains a finitely generated
$\fm_B$-primary ideal $J'\subset B$. Choose a finitely generated
$\fm_A$-primary ideal $J\subset A$ contained in the kernel of
the induced map $A\to B/J'$.

\begin{claim}\label{cl_reduce-to-bar}
$\lambda_B(M)=\lambda_{B/J'}(M)$ and
$\lambda_A(\psi_*M)=\lambda_{A/J}(\bar\psi_*M)$.
\end{claim}
\begin{pfclaim} Directly on the definition (and by applying
theorem \ref{th_about_K_0}(ii.b) to the surjection $B\to B/J'$)
we see that $\lambda_B(M)\leq\lambda_{B/J'}(M)$.
On the other hand, any surjection of $B$-modules $M'\to M$
with $M'$ finitely presented, factors through the natural
map $M'\to M'/J'M'$, and by theorem \ref{th_about_K_0}(i,ii.b)
we have $\lambda_B(M')-\lambda_{B/J'}(M'/J'M')=\lambda_B(J'M')\geq 0$,
whence the first stated identity. The proof of the second identity
is analogous.
\end{pfclaim}

In view of claim \ref{cl_reduce-to-bar} we are reduced to
proving the assertion for the morphism $\bar\psi$ and the
$B/J'$-module $M$, hence we may replace $\psi$ by $\bar\psi$,
and assume from start that $A$ and $B$ have Krull dimension zero.

\begin{claim}\label{cl_stanco}
Let $A$ be a measurable $K^+$-algebra of Krull dimension zero.
Then there exists a morphism $\phi:V\to A$ of measurable
$K^+$-algebras, with $V$ a valuation ring flat over $K^+$,
such that the residue field extension $\kappa(V)\to\kappa(A)$
is finite, and the induced map of value groups $\Gamma\to\Gamma_V$
is an isomorphism.
\end{claim}
\begin{pfclaim} Let $A_0\to A$ be a local ind-\'etale map,
from an object $A_0$ of $K^+\mAlg_0$. By lemma
\ref{lem_integral-fp-ext} we may find a $K^+$-flat valuation
ring $V_0$ in $K^+\mAlg_0$ and a local map $\phi_0:V\to A_0$,
inducing a finite residue field extension $\kappa(V_0)\to\kappa(A_0)$
and an isomorphism on value groups $\Gamma\isom\Gamma_V$.
Then $A_0$ has also Krull dimension zero; especially, it is
henselian, hence $\phi_0$ factors through a morphism
$\phi_0^\he:V_0^\he\to A_0$, from the henselization $V_0^\he$
of $V_0$. Let $E\subset\kappa(A)$ be the largest separable
subextension of $\kappa(V_0)=\kappa(V_0^\he)$ contained in
$\kappa(A)$; there exists a unique (up to unique isomorphism)
finite \'etale morphism $V_0^\he\to V'$ with an isomorphism
$\kappa(V')\isom E$ of $\kappa(V_0)$-algebras, and $\phi_0^\he$
factors through a morphism $V'\to A_0$.
Notice that $V'$ is still a henselian valuation ring and
a measurable $K^+$-algebra, hence we may replace $V_0$ by
$V'$, and assume that $V_0$ is henselian, and the residue
field extension $\bar\phi_0:\kappa(V_0)\to\kappa(A_0)$ is
purely inseparable. Since $A_0$ is henselian, we may write
$A$ as the colimit of a filtered system $(A_i~|~i\in I)$
of finite \'etale $A_0$-algebras. Now, on the one hand,
$\bar\phi_0$ induces an equivalence from the category
of finite \'etale $\kappa(V_0)$-algebras, to the category
of finite \'etale $\kappa(A_0)$-algebras (lemma
\ref{lem_replace}(i)). On the other hand, the category
of finite \'etale $V_0$-algebras is equivalent to the
category of finite \'etale $\kappa(V_0)$-algebras, and
likewise for $A_0$. Therefore, for every $i\in I$ we may
find a finite \'etale morphism $V_0\to V_i$, unique up
to unique isomorphism, inducing an isomorphism of
$\kappa(A_0)$-algebras :
\set\begin{equation}\label{eq_ouf}
\kappa(V_i)\otimes_{\kappa(V_0)}\kappa(A_0)\isom\kappa(A_i)
\end{equation}
and the transition maps of residue fields
$\kappa(A_i)\to\kappa(A_j)$ induce unique maps
$V_i\to V_j$ of $V_0$-algebras, compatible with the
isomorphisms \eqref{eq_ouf}. Hence, the resulting system
$(V_i~|~i\in I)$ is filtered, and its colimit is a valuation
ring $V$, which is still a measurable $K^+$-algebra.
Moreover, the field extensions $\kappa(V_i)\to\kappa(A_i)$
deduced from \eqref{eq_ouf} lift uniquely to maps of
$V_0$-algebras $V_i\to A_i$ (\cite[Ch.IV, Cor.18.5.12]{EGA4});
taking colimits, we get finally a map $V\to A$ as sought.
\end{pfclaim}

Let $\phi$ be as in claim \ref{cl_stanco}; clearly it
suffices to prove the sought identity for the two morphisms
$\psi\circ\phi$ and $\phi$, so we may replace $A$ by $V$,
and assume from start that $A$ is a valuation ring. Let us
set $d:=[\kappa(B):\kappa(A)]$; we deduce :
$$
\lambda_B(M)=
\inf\,\{d^{-1}\cdot\lambda_A(\psi_*M')~|~M'\in\cC_M\}\geq
d^{-1}\cdot\lambda_A(\psi_*M)
$$
by theorem \ref{th_about_K_0}(ii.b).
Furthermore, let us choose a surjection $B^{\oplus k}\to M$,
whose kernel we write in the form $K:=\bigcup_{i\in I}K_i$
where $(K_i~|~i\in I)$ is a filtered family of finitely
generated $B$-submodules of $K$. Next, by applying (i),
proposition \ref{prop_additive-vr}(i,iii.b), and theorem
\ref{th_about_K_0}(ii.b) we derive :
$$
\begin{aligned}
\lambda_B(M)\leq\, & \inf\,\{\lambda_B(B^{\oplus k}/K_i)~|~i\in I\}
=\inf\,\{d^{-1}\cdot\lambda_A(\psi_*(B^{\oplus k}/K_i))~|~i\in I\} \\
=\, & d^{-1}\cdot(\lambda_A(\psi_*B^{\oplus k})-
\sup\,\{\lambda_A(\psi_*K_i)~|~i\in I\}) \\
=\, & d^{-1}\cdot(\lambda_A(\psi_*B^{\oplus k})-\lambda_A(\psi_*K)) \\
=\, & d^{-1}\cdot\lambda_A(\psi_*M)
\end{aligned}
$$
whence the claim.

(iv): Let $(N_i~|~i\in I)$ be the filtered system of
finitely generated submodules of $M_2$. For every $i\in I$
we have short exact sequences :
$0\to M_1\cap N_i\to N_i\to \bar N_i\to 0$, where $\bar N_i$
is the image of $N_i$ in $M_3$. In view of the case of (ii)
already known, we may then replace $M_2$ by $N_i$, and thus
assume from start that $M_2$ is finitely generated, so that
we may find a finitely generated $\fm_A$-primary ideal
$J\subset A$ that annihilates $M_2$. By (iii) we have
$\lambda_A(M_2)=\lambda_{A/J}(M_2)$, and likewise for $M_1$
and $M_3$, hence we may replace $A$ by $A/J$. By claim
\ref{cl_stanco}, we may then find a morphism $V\to A$ of
measurable $K^+$-algebras with $V$ a valuation ring, inducing
a finite residue field extension $\kappa(V)\to\kappa(A)$;
then again (iii) reduces to the case where $A=V$, to which
one may apply (i) and proposition \ref{prop_additive-vr}(i)
to conclude the proof.

Next we consider assertion (ii) for the case where the
transition maps are surjective. We may assume that $I$ admits
a smallest element $i_0$; for every $i\in I$ let $K_i$
denote the kernel of the transition map $M_{i_0}\to M_i$.
We deduce a short exact sequence :
$$
0\to\bigcup_{i\in I}K_i\to M_{i_0}\to\colim_{i\in I}M_i\to 0
$$
and we may then compute using (iv) and the previous case of (ii) :
$$
\lambda_A(\colim_{i\in I}M_i)=
\lambda_A(M_{i_0})-\lambda_A(\bigcup_{i\in I}K_i)=
\lim_{i\in I}(\lambda_A(M_{i_0})-\lambda_A(K_i))=
\lim_{i\in I}\lambda_A(M_{i_0}/K_i)
$$
whence the claim.

(v): Since $\psi_*$ is an exact functor which commutes with
colimits, we may use (ii) to reduce to the case where $M$
is finitely presented, for which the assertion is already
known, in view of theorem \ref{th_about_K_0}(iv).
\end{proof}

\begin{remark}\label{rem_noeth-case}
Suppose that the valuation of $K$ is discrete; then one sees
easily that theorem \ref{th_about_K_0}(i,iii) still holds
(with simpler proof) when $A$ is replaced by any local noetherian
$K^+$-algebra, and by inspecting the definition, the resulting
map $\lambda_A$ is none else than the standard length function
for modules supported on $\{s(A)\}$.
\end{remark}

\begin{lemma}\label{lem_length-kills}
Let $A$ be a measurable $K^+$-algebra, and $M$ an $A$-module
supported at $s(A)$ with $\lambda_A(M)<\infty$. Then
$$
\{a\in K^+~|~\log|a|>\lambda_A(M)\}\subset\Ann_{K^+}M.
$$ 
\end{lemma}
\begin{proof} Using proposition \ref{prop_slight-gen}(ii),
we easily reduce, first, to the case where $M$ is finitely
generated, and second, to the case where $M$ is finitely
presented. Pick an ideal $I\subset A$ and a valuation ring
$V$ mapping onto $A/I$, as in lemma \ref{lem_find-val-ring};
by considering the $I$-adic filtration of $M$, the additivity
properties of $\lambda_A$ allow to further reduce to the
case where $M$ is an $A/I$-module. Next, by theorem
\ref{th_about_K_0}(ii.b) we may replace $A$ by $V$,
and therefore assume that $A$ is a valuation ring whose
value group equals $\Gamma$. In this case, $\lambda_A$
is computed by Fitting ideals, so the assertion follows
easily from \cite[Prop.6.3.6(iii)]{Ga-Ra}.
\end{proof}

\sset\subsubsection{}\label{subsec_last-step}
Let us consider a $K^+$-algebra $R_\infty$ that is the colimit
of an inductive system
\set\begin{equation}\label{eq_last-step}
R_0\to R_1\to R_2\to\cdots
\end{equation}
of morphisms of measurable $K^+$-algebras inducing integral ring
homomorphisms $\kappa(R_i)\to R_{i+1}/\fm_{R_i}R_i$ for every
$i\in\N$. The final step consists in generalizing the definition
of normalized length to the category $R_\infty\Mod_\bsb$ of
$R_\infty$-modules supported at the closed point $s(R_\infty)$
of $\Spec\,R_\infty$.
To this purpose, we shall axiomatize the general situation in
which we can solve this problem. Later we shall see that our axioms
are satisfied in many interesting cases.

\sset\subsubsection{}\label{subsec_axiom}
Hence, let $R_\infty$ be as in \eqref{subsec_last-step}.
After fixing an order-preserving isomorphism
\set\begin{equation}\label{eq_law-and-order}
(\Q\otimes_\Z\log\Gamma)^\wedge\isom\R
\end{equation}
we may regard the mappings $\lambda_A$ (for any measurable
$K^+$-algebra $A$) as real-valued functions on $A$-modules.
To ease notation, for every $R_n$-module $N$ supported on
$\{s(R_n)\}$ we shall write $\lambda_n(N)$ instead of
$\lambda_{R_n}(N)$.
Notice that, for every such $N$, and every $m\geq n$, the
$R_m$-module $R_m\otimes_{R_n}N$ is supported at $\{s(R_m)\}$,
since by assumption the map $R_n\to R_m/\fm_{R_n}R_m$ is
integral.

\begin{definition}\label{def_ind-measur}
In the situation of \eqref{subsec_axiom},
we say that $R_\infty$ is an {\em ind-measurable\/}
$K^+$-algebra, if there exists a sequence of real
{\em normalizing factors\/} $(d_n>0~|~n\in\N)$ such that:
\begin{enumerate}
\alphaenu
\item
For every $n\in\N$ and every object $N$ of $R_n\Mod_\cohs$,
the sequence :
$$
m\mapsto d_m^{-1}\cdot\lambda_m(R_m\otimes_{R_n}N)
$$
converges to an element
$\lambda_\infty(R_\infty\otimes_{R_n}N)\in\R$.
\item
For every $m\in\N$, every finitely generated $\fm_{R_0}$-primary
ideal $I\subset R_0$, and every $\eps>0$ there exists
$\delta(m,\eps,I)>0$ such that the following holds. For every
$n\in\N$ and every surjection $N\to N'$ of finitely presented
$R_n/IR_n$-modules generated by $m$ elements, such that
$$
|\lambda_\infty(R_\infty\otimes_{R_n}N)-
\lambda_\infty(R_\infty\otimes_{R_n}N')|\leq\delta(m,\eps,I)
$$
we have :
$$
d_n^{-1}\cdot|\lambda_n(N)-\lambda_n(N')|\leq\eps.
$$
\end{enumerate}
\end{definition}

\sset\subsubsection{}
Assume now that $R_\infty$ is ind-measurable, and let $N$ be
a finitely presented $R_n$-module supported on $\{s(R_n)\}$.
The first observation is that
$\lambda_\infty(R_\infty\otimes_{R_n}N)$ only depends on
the $R_\infty$-module $R_\infty\otimes_{R_n}N$. Indeed,
suppose that
$R_\infty\otimes_{R_m}M\simeq R_\infty\otimes_{R_n}N$
for some $m\in\N$ and some finitely presented $R_m$-module
$M$; then there exists $p\geq m,n$ such that
$R_p\otimes_{R_m}M\simeq R_p\otimes_{R_n}N$,
and then the assertion is clear.

The second observation -- contained in the following lemma
\ref{lem_restriction} -- will show that the conditions of
definition \ref{def_ind-measur} impose some non-trivial
restrictions on the inductive system $(R_n~|~n\in\N)$.

\begin{lemma}\label{lem_restriction}
Let $(R_n; d_n~|~n\in\N)$ be an inductive system of
measurable $K^+$-algebras and a sequence of positive
reals, fulfilling conditions {\em (a)} and {\em (b)}
of definition {\em\ref{def_ind-measur}}. Then:
\begin{enumerate}
\item
The natural map
$$
M\to R_m\otimes_{R_n}M
$$
is injective for every $n,m\in\N$ with $m\geq n$ and every
$R_n$-module $M$.
\item
Especially, the transition maps $R_n\to R_{n+1}$ are injective
for every $n\in\N$.
\item
Suppose that $(d'_n~|~n\in\N)$ is another sequence of
positive reals such that conditions {\em (a)} and {\em (b)}
hold for the datum $(R_n; d'_n~|~n\in\N)$.
Then the sequence $(d_n/d'_n~|~n\in\N)$ converges to
a non-zero real number.
\end{enumerate}
\end{lemma}
\begin{proof} (i): We reduce easily to the case where $M$ is
finitely presented over $R_n$. Let
$N\subset\Ker\,(M\to R_m\otimes_{R_n}M)$ be a finitely generated
$R_n$-module; since $R_n$ is coherent, $N$ is a finitely presented
$R_n$-module. We suppose first that $M$ is in $R_n\Mod_\cohs$.

\begin{claim}\label{cl_represent} The natural map
$R_m\otimes_{R_n}M\to R_m\otimes_{R_n}(M/N)$ is an isomorphism.
\end{claim}
\begin{pfclaim} On the one hand, the $R_m$-module $R_m\otimes_{R_n}M$
represents the functor
$$
R_m\Mod\to\Set\quad:\quad Q\mapsto\Hom_{R_n}(M,Q).
$$
On the other hand, the assumption on $N$ implies that
$\Hom_{R_n}(M,Q)=\Hom_{R_n}(M/N,Q)$ for every $R_m$-module
$Q$, so the claim follows easily.
\end{pfclaim}

From claim \ref{cl_represent} we deduce that the natural
map $R_\infty\otimes_{R_n}M\to R_\infty\otimes_{R_n}(M/N)$
is an isomorphism, and then condition (b) says that
$\lambda_n(M)=\lambda_n(M/N)$, hence $\lambda_n(N)=0$ and
finally $N=0$, as stated. Next, suppose $M$ is any finitely
presented $R_n$-module, and pick a finitely generated
$\fm_{R_n}$-primary ideal $I\subset R_n$.

\begin{claim}\label{cl_susan}
There exists $c\in\N$ such that $N\cap I^{k+c}M=I^k(N\cap I^cM)$
for every $k\geq 0$.
\end{claim}
\begin{pfclaim} We may find a local ind-\'etale map
$A\to R_n$ of $K^+$-algebras, where $A$ is an object
of $K^+\mAlg_0$, and such that $I$, $M$ and
$N$ descend respectively to a finitely generated ideal
$I_0\subset A$, a finitely presented $A$-module $M_0$,
and a finitely generated submodule $N_0\subset M_0$.
Then theorem \ref{th_Rees} ensures the existence of
$c\in\N$ such that
$N_0\cap I_0^{k+c}M_0=I_0^k(N_0\cap I_0^cM_0)$ for
every $k\geq 0$. Since $R_n$ is a faithfully flat
$A$-algebra, the claim follows.
\end{pfclaim}

Pick $c\in\N$ as in claim \ref{cl_susan}; then
$N/(N\cap I^{k+c}M)$ is in the kernel of the natural map
$M/I^{k+c}M\to R_m\otimes_{R_n}(M/I^{k+c}M)$, hence
$N=N\cap I^{k+c}M$ by the foregoing, so that $N\subset I^kN$
for every $k\geq 0$, and finally $N=0$ by Nakayama's lemma.

(ii) is a special case of (i). To show (iii), let us denote
by $\lambda'_\infty(M)$ the normalized length of any object
$M$ of $R_\infty\Mod_\cohs$, defined using the sequence
$(d'_n~|~n\in\N)$. From (b) it is clear that
$\lambda_\infty(M),\lambda'_\infty(M)\neq 0$ whenever $M\neq 0$.
Then, for any such non-zero $M$, the quotient
$\lambda'_\infty(M)/\lambda_\infty(M)$ is the limit of the
sequence $(d_n/d'_n~|~n\in\N)$.
\end{proof}

\begin{remark}\label{rem_ind-measure}
(i)\ \ 
There is another situation of interest which leads to a
well-behaved notion of normalized length. Namely, suppose
that $R_\bullet:=(R_n~|~n\in\N)$ is an inductive system of
local homomorphisms of local noetherian rings, such that the
fibres of the induced morphisms
$\Spec\,R_{n+1}\to\Spec\,R_n$ have dimension zero, and denote
by $R_\infty$ the inductive limit of the system $R_\bullet$.
For every $n\in\N$, let also $\lambda_n$ be the usual length
function on the set of isomorphism classes of finitely
generated $R_n$-modules supported on $s(R_n)$.
Then for every $m,n\in\N$ with $m\geq n$, and every $R_n$-module
$M$ of finite length, the $R_m$-module $R_m\otimes_{R_n}M$
has again finite length, so the analogues of conditions
(a) and (b) of \eqref{subsec_axiom} can be formulated
(cp. remark \ref{rem_noeth-case}), and if these conditions
hold for $R_\bullet$, we shall say that $R_\infty$ is an
{\em ind-measurable ring}. In such situation, lemma
\ref{lem_restriction} -- as well as the forthcoming lemma
\ref{lem_two-defs-coincide} and theorem \ref{th_most-general} --
still hold, with simpler proofs : we leave the details to the
reader.

(ii)\ \
In spite of the uniqueness properties expressed by lemma
\ref{lem_restriction}, we do not know to which extent the
normalized length of an ind-measurable $K^+$-algebra
depends on the chosen tower of measurable algebras.
Namely, suppose that $(R_n~|~n\in\N)$ and $(R'_n~|~n\in\N)$
are two such towers, with isomorphic colimit $R_\infty$,
and suppose that we have found normalizing factors
$(d_n~|~n\in\N)$ (resp. $(d'_n~|~n\in\N)$) for the
first (resp. second) tower, whence a normalized length
$\lambda_\infty$ (resp. $\lambda'_\infty$) for
$R_\infty$-modules. Then we do not know whether the
ratio of $\lambda_\infty$ and $\lambda'_\infty$ is
a constant.
\end{remark}

Next, for a given finitely generated $R_\infty$-module $M$
supported on $\{s(R_\infty)\}$, we shall proceed as in
\eqref{subsec_C_M} : we denote by $\cC_M$ the set of isomorphism
classes of finitely presented $R_\infty$-modules supported on
$\{s(R_\infty)\}$ that admit a surjection $M'\to M$, and we set
$$
\lambda^*_\infty(M):=\inf\,\{\lambda_\infty(M')~|~M'\in\cC_M\}
\in\R.
$$
Directly on the definitions, one checks that
$\lambda^*_\infty(M)=\lambda_\infty(M)$ if $M$ is finitely presented.

\begin{lemma}\label{lem_two-defs-coincide}
Let $M$ be a finitely generated $R_\infty$-module,
$\Sigma\subset M$ any finite set of generators, $(N_i~|~i\in I)$
any filtered system of objects of $R_\infty\Mod_\bsb$, with surjective
transition maps, such that\/ $\colim_{i\in I}N_i\simeq M$. Then :
\begin{enumerate}
\item
$\lambda^*_\infty(M)=
\displaystyle{\lim_{n\to\infty}}
d_n^{-1}\cdot\lambda_n(\Sigma R_n)$.
\item
If every $N_i$ is finitely generated, we have :
$$
\inf\,\{\lambda^*_\infty(N_i)~|~i\in I\}=\lambda^*_\infty(M).
$$
\end{enumerate}
\end{lemma}
\begin{proof} To start out, we show assertion (ii) in the
special case where all the modules $N_i$ are finitely
presented. Indeed, let us pick any surjection
$\phi:M'\to M$ with $M'\in\cC_M$. We may find $i\in I$
such that $\phi$ lifts to a map $\phi_i:M'\to N_i$
(\cite[Prop.2.3.16(ii)]{Ga-Ra}), and up to replacing $I$
by a cofinal subset, we may assume that $\phi_i$ is defined
for every $i\in I$. Moreover
$$
\colim_{i\in I}\Coker\,\phi_i=\Coker\,\phi=0
$$
hence there exists $i\in I$ such that $\phi_i$ is surjective.
It follows easily that $\lambda_\infty(N_i)\leq\lambda_\infty(M')$,
whence (ii) in this case.

(i): Let $k$ be the cardinality of $\Sigma$, $\eps>0$ any
real number, $Q\subset\Ann_{R_0}\Sigma$ a finitely generated
$\fm_{R_0}$-primary ideal, $S:=R_\infty/QR_\infty$ and
$\beta:S^{\oplus k}\to M$ a surjection that sends the standard
basis onto $\Sigma$. By filtering $\Ker\,\beta$ by the system
$(K_j~|~j\in J)$ of its finitely generated submodules, we
obtain a filtered system $(N_j:=S^{\oplus k}/K_j~|~j\in J)$
of finitely presented $R_\infty$-modules supported on $\{s(R_\infty)\}$,
with surjective transition maps and colimit isomorphic to $M$.
Let $\psi_j:N_j\to M$ be the natural map; by the foregoing,
we may find $j_0\in J$ such that
\set\begin{equation}\label{eq_good-approx}
0\leq\lambda_\infty(N_{j_0})-\lambda^*_\infty(M)
\leq\min(\delta(k,\eps,Q),\eps)
\end{equation}
(notation of definition \ref{def_ind-measur}(b)).
We can then find $n\in\N$ and a finitely presented $R_n$-module
$M'_n$ supported on $\{s(R_n)\}$ and generated by at most $k$
elements, such that $N_{j_0}=R_\infty\otimes_{R_n}M'_n$; we set
$M'_m:=R_m\otimes_{R_n}M_n$ for every $m\geq n$ and let $M_m$ be
the image of $M'_m$ in $M$.
Up to replacing $n$ by a larger integer, we may assume that
$M_m=\Sigma R_m$ for every $m\geq n$.
Choose $m\in\N$ so that :
\set\begin{equation}\label{eq_first_eps}
|\lambda_\infty(N_{j_0})-d^{-1}_m\cdot\lambda_m(M'_m)|\leq\eps.
\end{equation}
Let $(M_{m,i}~|~i\in I)$ be a filtered system of finitely presented
$R_m$-modules supported at $\{s(R_m)\}$, with surjective transition
maps, such that $\colim_{i\in I}M_{m,i}=M_m$. Arguing as in the foregoing,
we show that there exists $i\in I$ such that the natural map
$\phi:M'_m\to M_m$ factors through a surjection
$\phi_i:M'_m\to M_{m,i}$, and up to replacing $I$ by a
cofinal subset, we may assume that such a surjection $\phi_i$
exists for every $i\in I$. We obtain therefore a compatible
system of surjections of $R_\infty$-modules :
$$
N_{j_0}\simeq R_\infty\otimes_{R_m}M'_m\to
R_\infty\otimes_{R_m}M_{m,i}\to M
$$
and combining with \eqref{eq_good-approx} we find :
$$
|\lambda_\infty(N_{j_0})-\lambda_\infty(R_\infty\otimes_{R_m}M_{m,i})|
\leq\delta(k,\eps,Q)
\qquad\text{for every $i\in I$}.
$$
In such situation, condition (b) of definition \ref{def_ind-measur}
ensures that :
$$
d_m^{-1}\cdot|\lambda_m(M'_m)-\lambda_m(M_{m,i})|\leq\eps
\qquad\text{for every $i\in I$}.
$$
Therefore : $d_m^{-1}\cdot|\lambda_m(M'_m)-\lambda_m(M_m)|\leq\eps$,
by proposition \ref{prop_slight-gen}(ii). Combining with
\eqref{eq_first_eps} and again \eqref{eq_good-approx} we obtain :
$$
|\lambda^*_\infty(M)-d_m^{-1}\cdot\lambda_m(M_m)|\leq 3\eps
$$
which implies (i).

(ii): With no loss of generality, we may assume that $I$ admits
a smallest element $i_0$. Let us fix a surjection
$F:=R_\infty^{\oplus k}\to N_{i_0}$, and for every $i\in I$, let
$C_i$ denote the kernel of the induced surjection $F\to N_i$.
We consider the filtered system $(D_j~|~j\in J)$ consisting of all
finitely generated submodules $D_j\subset F$ such that $D_j\subset C_i$
for some $i\in I$. It is clear that $\colim_{j\in J}F/D_j\simeq M$,
hence $\lambda^*_\infty(M)=\inf\,\{\lambda_\infty(F/D_j)~|~j\in J\}$,
by (i). On the other, by construction, for every $j\in J$ we may
find $i\in I$ such that $\lambda^*_\infty(N_i)\leq\lambda_\infty(F/D_j)$;
since clearly $\lambda^*_\infty(N_i)\geq\lambda^*_\infty(M)$ for
every $i\in I$, the assertion follows.
\end{proof}

\sset\subsubsection{}
Let now $M$ be an arbitrary object of the category
$R_\infty\Mod_\bsb$. We let :
$$
\lambda_\infty(M):=\sup\,\{\lambda^*_\infty(M')~|~M'\subset M
\text{ and $M'$ is finitely generated}\}\in\R\cup\{+\infty\}.
$$

\begin{lemma}\label{lem_agree-previous}
If $M$ is finitely generated, then
$\lambda_\infty(M)=\lambda^*_\infty(M)$.
\end{lemma}
\begin{proof} Let $N\subset M$ be any finitely generated
submodule. We choose finite sets of generators $\Sigma\subset N$
and $\Sigma'\subset M$ with $\Sigma\subset\Sigma'$. In view
of proposition \ref{prop_slight-gen}(iv) we have
$\lambda_n(\Sigma R_n)\leq\lambda_n(\Sigma'R_n)$ for
every $n\in\N$, hence $\lambda^*_\infty(N)\leq\lambda^*_\infty(M)$,
by lemma \ref{lem_two-defs-coincide}(i). The contention follows
easily.
\end{proof}

\begin{theorem}\label{th_most-general}
{\em (i)}\ \ Let $(M_i~|~i\in I)$ be a filtered
system of objects of $R_\infty\Mod_\bsb$, and suppose that either :
\begin{enumerate}
\alphaenu
\item
all the transition maps of the system are injections, or
\item
all the transition maps are surjections and
$\lambda_\infty(M_i)<+\infty$ for every $i\in I$.
\end{enumerate}
Then :
$$
\lambda_\infty(\colim_{i\in I}M_i)=\lim_{i\in I}\lambda_\infty(M_i).
$$
\begin{enumerate}
\addenu
\item
Let $0\to M'\to M\to M''\to 0$ be a short exact sequence in
$R_\infty\Mod_\bsb$. Then :
$$
\lambda_\infty(M)=\lambda_\infty(M')+\lambda_\infty(M'').
$$
\item
Let $M$ be a finitely presented $R_\infty$-module,
$N\subset M$ a submodule supported at $s(R_\infty)$.
Then $\lambda_\infty(N)=0$ if and only if $N=0$.
\end{enumerate}
\end{theorem}
\begin{proof} The proof of (i) in case (a) is the same as that of
proposition \ref{prop_slight-gen}(iii.b).

(ii): We proceed in several steps :

$\bullet$\ \ 
Suppose first that $M$ and $M''$ are finitely presented,
hence $M'$ is finitely generated (\cite[Lemma 2.3.18(ii)]{Ga-Ra}).
We may then find an integer $n\in\N$ and finitely presented
$R_n$-modules $M_n$ and $M''_n$ such that
$M\simeq R_\infty\otimes_{R_n}M_n$ and
$M''\simeq R_\infty\otimes_{R_n}M''_n$. For every $m\geq n$ we
set $M_m:=R_m\otimes_{R_n}M_n$ and likewise we define $M''_m$.
Up to replacing $n$ by a larger integer, we may assume that
the given map $\phi:M\to M''$ descends to a surjection
$\phi_n:M_n\to M''_n$, and then $\phi$ is the colimit of the
induced maps $\phi_m:=\one_{R_m}\otimes_{R_n}\phi_n$, for every
$m\geq n$. Moreover, $M'\simeq\colim_{m\geq n}\Ker\,\phi_m$,
and by the right exactness of the tensor product, the image of
$\Ker\,\phi_m$ generates $M'$ for every $m\geq n$.
Furthermore, since the natural maps $M_m\to M_p$ are injective
for every $p\geq m\geq n$ (lemma \ref{lem_restriction}), the
same holds for the induced maps $\Ker\,\phi_m\to\Ker\,\phi_p$.
The latter factors as a composition :
$$
\Ker\,\phi_m\stackrel{\alpha}{\longrightarrow}
R_p\otimes_{R_m}\Ker\,\phi_m\stackrel{\beta}{\longrightarrow}
\Ker\,\phi_p
$$
where $\alpha$ is injective (lemma \ref{lem_restriction})
and $\beta$ is surjective. In other words,
$R_p\cdot\Ker\,\phi_m=\Ker\,\phi_p$ for every $p\geq m\geq n$.
In such situation, lemma \ref{lem_two-defs-coincide}(i)
ensures that :
$$
\lambda_\infty(M')=
\lim_{m\geq n}d^{-1}_m\cdot\lambda_m(\Ker\,\phi_m)
$$
and likewise :
$$
\lambda_\infty(M)=
\lim_{m\geq n}d^{-1}_m\cdot\lambda_m(M_m)\qquad
\lambda_\infty(M'')=
\lim_{m\geq n}d^{-1}_m\cdot\lambda_m(M''_m).
$$
To conclude the proof of (ii) in this case, it suffices
then to apply proposition \ref{prop_slight-gen}(iv).

$\bullet$\ \ 
Suppose next that $M$, $M'$ (and hence $M''$) are finitely
generated. We choose a filtered system $(M_i~|~i\in I)$
of finitely presented $R_\infty$-modules, with surjective
transition maps, such that $M\simeq\colim_{i\in I}M_i$.
After replacing $I$ by a cofinal subset, we may assume
that $M'$ is generated by a finitely generated submodule
$M'_i$ of $M_i$, for every $i\in I$, and that $(M'_i~|~i\in I)$
forms a filtered system with surjective transition maps,
whose colimit is necessarily $M'$; set also $M''_i:=M_i/M'_i$
for every $i\in I$, so that the colimit of the filtered
system $(M''_i~|~i\in I)$ is $M''$. In view of lemmata
\ref{lem_two-defs-coincide}(ii) and \ref{lem_agree-previous},
we are then reduced to showing the identity :
$$
\inf\,\{\lambda_\infty^*(M_i)~|~i\in I\}=
\inf\,\{\lambda_\infty^*(M'_i)~|~i\in I\}+
\inf\,\{\lambda_\infty^*(M''_i)~|~i\in I\}
$$
which follows easily from the previous case.

$\bullet$ Suppose now that $M$ is finitely generated.
We let $(M'_i~|~i\in I)$ be the filtered family of finitely
generated submodules of $M'$. Then :
$$
M'\simeq\colim_{i\in I} M'_i \quad\text{and}\quad
M''\simeq\colim_{i\in I} M/M'_i.
$$
Hence :
$$
\lambda_\infty(M')=\lim_{i\in I}\lambda_\infty(M'_i)
\quad\text{(resp.\quad
$\lambda_\infty(M'')=\lim_{i\in I}\lambda_\infty(M/M'_i)$)}
$$
by (i.a) (resp. by lemmata \ref{lem_two-defs-coincide}(ii) and
\ref{lem_agree-previous}). However, the foregoing case shows that
$\lambda_\infty(M)=\lambda_\infty(M'_i)+\lambda_\infty(M/M'_i)$
for every $i\in I$, so assertion (ii) holds also in this case.

$\bullet$ Finally we deal with the general case. Let
$(M_i~|~i\in I)$ be the filtered system of finitely
generated submodules of $M$; we denote by $M''_i$ the
image of $M_i$ in $M''$, and set $M'_i:=M'\cap M_i$
for every $i\in I$. By (i.a) we have :
$$
\lambda_\infty(M)=\lim_{i\in I}\lambda_\infty(M_i)
$$
and likewise for $M'$ and $M''$. Since we already know
that $\lambda_\infty(M_i)=
\lambda_\infty(M'_i)+\lambda_\infty(M''_i)$ for every
$i\in I$, we are done.

(iii): Let $f\in N$ be any element; in view of (ii) we
see that $\lambda_\infty(fR_\infty)=0$, and it suffices
to show that $f=0$. However, we may find $n\in\N$
and a finitely presented $R_n$-module $M_n$ such that
$M\simeq R_\infty\otimes_{R_n}M_n$; notice that the
natural map $M_n\to M$ is injective, by lemma
\ref{lem_restriction}(i). We may also assume that $f$
is in the image of $M_n$. Let $I\subset\Ann_{R_n}(f)$
be a finitely generated $\fm_{R_n}$-primary ideal; after
replacing $M$ by $M/IM$, we may assume that $M_n$ is
supported at $s(R_n)$. In light of (ii), we see
that $\lambda_\infty(M)=\lambda_\infty(M/fR_\infty)$,
hence $\lambda_n(M_n)=\lambda_n(M_n/fR_n)$, due to
condition (b) of definition \ref{def_ind-measur}. Hence
$\lambda_n(fR_n)=0$ by proposition \ref{prop_slight-gen}(iv),
and finally $f=0$ by theorem \ref{th_about_K_0}(i,iii.a).

To conclude, we consider assertion (i) in case (b) : set
$M:=\colim_{i\in I}M_i$; it is clear that
$\lambda_\infty(M)\leq\lambda_\infty(M_i)\leq\lambda_\infty(M_j)$
whenever $i\geq j$, hence
$$
\lim_{i\in I}\lambda_\infty(M_i)=\inf\,\{\lambda_\infty(M_i)~|~i\in I\}
\geq\lambda_\infty(M).
$$
For the converse inequality, fix $\eps>0$; without loss
of generality, we may assume that $I$ admits a smallest
element $i_0$, and we can find a finitely generated submodule
$N_{i_0}\subset M_{i_0}$ such that
$\lambda_\infty(M_{i_0})-\lambda_\infty(N_{i_0})<\eps$.
For every $i\in I$, let $N_i\subset M_i$ be the image of $N_{i_0}$,
and let $N\subset M$ be the colimit of the filtered system
$(N_i~|~i\in I)$; then $M_i/N_i$ is a quotient of $M_{i_0}/N_{i_0}$,
and the additivity assertion (ii) implies that
$\lambda_\infty(M_i)-\lambda_\infty(N_i)<\eps$ for every
$i\in I$.
According to lemma \ref{lem_two-defs-coincide}(ii) (and lemma
\ref{lem_agree-previous}) we have :
$$
\lambda_\infty(M)\geq\lambda_\infty(N)=
\inf\,\{\lambda_\infty(N_i)~|~i\in I\}\geq
\inf\,\{\lambda_\infty(M_i)~|~i\in I\}-\eps
$$
whence the claim.
\end{proof}

\sset\subsubsection{}
We wish now to show that the definition of normalized length
descends to almost modules (see \eqref{sec_norm-lengths}).
Namely, we have the following :

\begin{proposition}\label{prop_desc-to-alm}
Let $M$, $N$ be two objects of $R_\infty\Mod_\bsb$ such that
$M^a\simeq N^a$. Then $\lambda_\infty(M)=\lambda_\infty(N)$.
\end{proposition}
\begin{proof} Using additivity (theorem \ref{th_most-general}(ii)),
we easily reduce to the case where $M^a=0$, in which case
we need to show that $\lambda_\infty(M)=0$. Using theorem
\ref{th_most-general}(i) we may further assume that $M$
is finitely generated. Then, in view of lemma
\ref{lem_two-defs-coincide}(i), we are reduced to showing
the following :

\begin{claim} Let $A$ be any measurable $K^+$-algebra,
and $M$ any object of $A\Mod_\bsb$ such that $M^a=0$.
Then $\lambda_A(M)=0$.
\end{claim}
\begin{pfclaim}[] Arguing as in the foregoing we reduce to
the case where $M$ is finitely generated. Then, let
us pick $I\subset A$ and $\phi:V\to A/I$ as in lemma
\ref{lem_find-val-ring}. It suffices to show the assertion
for the finitely generated module $\bigoplus_{n\in\N}I^nM/I^{n+1}M$,
hence we may assume that $I\subset\Ann_AM$, in which case,
by proposition \ref{prop_slight-gen}(iii) we may replace
$A$ by $V$ and assume throughout that $A$ is a valuation
ring. Then the claim follows from propositions \ref{prop_slight-gen}(i)
and \ref{prop_additive-vr}(ii).
\end{pfclaim}
\end{proof}

\sset\subsubsection{}\label{subsec_extend-to-alm}
Proposition \ref{prop_desc-to-alm} suggests the following
definition. We let $R^a_\infty\Mod_\bsb$ be the full
subcategory of $R^a_\infty\Mod$ consisting of all the
$R^a_\infty$-modules $M$ such that $M_!$ is supported
at $s(R_\infty)$, in which case we say that $M$ {\em is
supported at\/} $s(R_\infty)$. Then, for every such $M$ we set :
$$
\lambda_\infty(M):=\lambda_\infty(M_!).
$$
With this definition, it is clear that theorem
\ref{th_most-general}(i,ii) extends {\em mutatis mutandis\/}
to almost modules. For future reference we point out :

\begin{lemma}\label{lem_fut-refer}
Let $0\to M'\to M\to M''\to 0$ be a short exact sequence of
$R^a_\infty$-modules. We have :
\begin{enumerate}
\item
If $M$ lies in $R^a_\infty\Mod_\bsb$, then \ \
$\lambda_\infty(abM)\leq\lambda_\infty(aM')+\lambda_\infty(bM'')$
\ \ for every $a,b\in\fm$.
\item
If $M$ is almost finitely presented, and $M'$ is supported at
$s(R_\infty)$, then $\lambda_\infty(M')=0$
if and only if $M'=0$.
\end{enumerate}
\end{lemma}
\begin{proof} (i): To start out, we deduce a short exact sequence :
$0\to bM\cap M'\to bM\to bM''\to 0$.
Next, let $N:=\Ker(a:bM\to bM)$, and denote by $N'$ the
image of $N$ in $bM''$; there follows a short exact sequence :
$0\to a(bM\cap M')\to abM\to bM''/N'\to 0$. The claim follows.

(ii): We reduce easily to the case where $M'$ is a cyclic
$R^a_\infty$-module, say $M'=R_\infty^ax$ for some $x\in M_!$.
In this case, let $I\subset\Ann_{R_0}(x)$ be a finitely
generated $\fm_{R_0}$-primary ideal; we may then replace
$M$ by $M/IM$, and assume that $M$ lies in $R^a_\infty\Mod_\bsb$
as well. We remark :

\begin{claim}\label{cl_counter-dir}
$M$ is almost finitely presented if and only if, for every
$b\in\fm$ there exist a finitely presented $R_\infty$-module
$N$ and a morphism $M\to N^a$ whose kernel and cokernel are
annihilated by $b$. Moreover if $M$ is supported at $s(R_\infty)$,
one can choose $N$ to be supported at $s(R_\infty)$.
\end{claim}
\begin{pfclaim} The ``if'' direction is clear. For the
``only if'' part, we use \cite[Cor.2.3.13]{Ga-Ra}, which
provides us with a morphism $\phi:N\to M_!$ with $N$ finitely
presented over $R_\infty$, such that
$b\cdot\Ker\,\phi=b\cdot\Coker\,\phi=0$. Then $b\cdot\one_N^a$
factors through a morphism $\phi':\Img\,\phi^a\to N^a$, and
$b\cdot\one_M$ factors through a morphism $\phi'':M\to\Img\,\phi^a$;
the kernel and cokernel of $\phi'\circ\phi''$ are annihilated
by $b^2$. Finally, suppose that $M$ is supported on $s(R_\infty)$,
and let $I\subset R_\infty$ be any finitely generated ideal
such that $V(I)=\{s(R_\infty)\}$. By assumption, for every
$f\in I$ and every $m\in M_!$ there exists $n\in\N$ such that
$f^nm=0$; it follows that $I^n$ annihilates $\Img\,\phi$ for
every sufficiently large $n\in\N$, and we may then replace $N$
by $N/I^nN$.
\end{pfclaim}

Let $M$ and $M'$ be as in (ii), and choose a morphism $\phi:M\to N^a$
as in claim \ref{cl_counter-dir}; by adjunction we get a map
$\psi:M'_!\to M_!\to N$ with $b\cdot\fm\cdot\Ker\,\psi=0$. It follows
that $\lambda_\infty(\Img\,\psi)=0$, hence $\Img\,\psi=0$, by theorem
\ref{th_most-general}(iii). Hence $bM'=0$; since $b$ is arbitrary,
the assertion follows.
\end{proof}

Simple examples show that an almost finitely generated
(or even almost finitely presented) $R^a_\infty$-module
may fail to have finite normalized length. The useful finiteness
condition for almost modules is contained in the following :

\begin{definition}\label{def_alm-finite-length}
Let $M$ be a $R^a_\infty$-module supported at $s(R_\infty)$.
We say that $M$ has {\em almost finite length\/} if
$\lambda_\infty(bM)<+\infty$ for every $b\in\fm$.
\end{definition}

\begin{lemma}\label{lem_alm-finite-lambda}
{\em (i)}\ \ The set of isomorphism classes of $R^a_\infty$-modules
of almost finite length forms a closed subset of the uniform space
$\cM(A)$ (notation of \cite[\S2.3]{Ga-Ra}).
\begin{enumerate}
\addenu
\item
Especially, every almost finitely generated $R^a_\infty$-module
supported at $s(R_\infty)$ has almost finite length.
\end{enumerate}
\end{lemma}
\begin{proof} Assertion (i) boils down to the following :

\begin{claim} Let $a,b\in\fm$, and $f:N\to M$, $g:N\to M'$
morphisms of $R^a_\infty$-modules, such that the kernel and
cokernel of $f$ and $g$ are annihilated by $a\in\fm$, and
such that $\lambda_\infty(bM)<+\infty$. Then
$\lambda_\infty(a^2bM')<+\infty$.
\end{claim}
\begin{pfclaim} By assumption, $\Ker\,f\subset\Ker\,a\cdot\one_N$,
whence an epimorphism $f(N)\to aN$; likewise, we have an epimorphism
$g(N)\to aM'$. It follows that
$$
\lambda_\infty(a^2bM')\leq\lambda_\infty(ab\cdot g(M))\leq
\lambda_\infty(abN)\leq\lambda_\infty(b\cdot f(N))<+\infty
$$
as stated.
\end{pfclaim}

(ii) follows from (i) and the obvious fact that every finitely
generated $R_\infty$-module supported at $s(R_\infty)$ has finite
normalized length.
\end{proof}

\sset\subsubsection{}\label{subsec_uniform}
For the case of a measurable $K^+$-algebra $A$ of dimension
zero, we can show a further {\em Lipschitz type\/} uniform
estimate for the normalized length. Namely, for any integer
$k>0$ define the set $\cM_k(A^a)$ with its uniform structure
as in \eqref{subsec_general-uniformity}; {\em i.e.} we have
the fundamental system of entourages
$$
(E_r~|~r\in\R_{>0})
$$
where each $E_r$ consists of the pairs $(N,N')$ such
that there exist a third $A^a$-module $N''$ and $A^a$-linear
maps $N''\to N$, $N''\to N'$ whose kernel and cokernels are
annihilated by any $b\in K^+$ such that $\log|b|\geq r$.

\begin{lemma}\label{lem_ta-ta}
In the situation of \eqref{subsec_uniform}, we have :
$$
|\lambda_{A^a}(N)-\lambda_{A^a}(N')|\leq
4k\cdot\length_A(A/\fm A)\cdot r
$$
for every $r\in\R_{>0}$ and every $(N,N')\in E_r$.
\end{lemma}
\begin{proof} We begin with the following

\begin{claim}\label{cl_lipstick}
For any finitely generated $A$-module $N$, and every
$b\in K^+\setminus\{0\}$ we have
$$
\lambda_A(N/bN)\leq\length_A(N/\fm N)\cdot\log|b|.
$$
\end{claim}
\begin{pfclaim} Choose a map $\psi:V\to A$ of measurable
$K^+$-algebras, inducing a finite residue field extension
$\kappa(V)\to\kappa(A)$, where $V$ is a $K^+$-flat valuation
ring, and the induced map of value groups $\Gamma\to\Gamma_V$
is an isomorphism (claim \ref{cl_stanco}).
Let $d:=\dim_{\kappa(V)}N/\fm_VN$; applying Nakayama's
lemma, we get a surjection $(V/bV)^{\oplus d}\to N/bN$
of $V$-modules, for every $b\in K^+$. Hence
$$
\lambda_A(N/bN)=
\frac{\lambda_V(\psi_*(N/bN))}{[\kappa(V):\kappa(A)]}\leq
\frac{d\log|b|}{[\kappa(V):\kappa(A)]}=
\length_A(N/\fm N)\cdot\log|b|
$$
as stated.
\end{pfclaim}

Now, suppose that $(N,N')\in E_{\log|b|}$ for some $b\in K^+$,
and pick maps $\phi:N''\to N$, $\psi:N''\to N'$ whose kernel
and cokernel are annihilated by $b$. Especially, if
$n_1,\dots,n_k\in N_*$ (resp. $n'_1,\dots,n'_k\in N'_*$)
is a system of generators for $N$ (resp. for $N'$), we may
find $n''_1,\dots,n''_{2k}\in N''_*$ such that
$\phi(n''_i)=b^2n_i$ for $i=1,\dots,k$ and $\psi(n''_i)=b^2n'_{i-k}$
for $i=k+1,\dots,2k$. After replacing $N''$ by its submodule
generated by $n''_1,\dots,n''_{2k}$, we may assume that
$N''\in\cM_{2k}(A)$, $\Ker\,\phi$ is still annihilated by
$b$, but $\Coker\,\phi$ is only annihilated by $b^2$. On
the one hand, we deduce that
$$
\lambda_{A^a}(N'')\geq\lambda_{A^a}(\phi(N''))\geq
\lambda_{A^a}(b^2N)=\lambda_{A^a}(N)-\lambda_{A^a}(N/b^2N)
$$
therefore
$$
\lambda_{A^a}(N)-\lambda_{A^a}(N'')\leq\lambda_{A^a}(N/b^2N).
$$
On the other hand, the map $N''\to N''$ given by the rule :
$n''\mapsto bn''$ for every $n''\in N''_*$, factors through
$\phi(N'')$, therefore $\lambda_{A^a}(N)\geq\lambda_{A^a}(bN'')$
so the same calculation yields the inequality
$$
\lambda_{A^a}(N'')-\lambda_{A^a}(N)\leq\lambda_{A^a}(N''/bN'').
$$
Taking into account claim \ref{cl_lipstick} (and proposition
\ref{prop_desc-to-alm}) we see that
$$
|\lambda_{A^a}(N'')-\lambda_{A^a}(N)|\leq
2k\cdot\length_A(A/\fm A)\cdot\log|b|.
$$
Of course, the same holds with $N$ replaced by $N'$, and
the lemma follows.
\end{proof}

We conclude this section with some basic examples of
the situation contemplated in definition \ref{def_ind-measur}.

\begin{example}\label{ex_flat-case}
Suppose that all the transition maps $R_n\to R_{n+1}$
of the inductive system in \eqref{subsec_axiom} are flat.
Then, proposition \ref{prop_slight-gen}(v) implies that
conditions (a) and (b) hold with:
$$
d_n:=\length_{R_n}(R_n/\fm_{R_0}R_n)
\qquad\text{for every $n\in\N$}.
$$
\end{example}

\begin{example}\label{ex_smooth}
Suppose that $p:=\chara\,\kappa>0$. Let $d\in\N$ be any integer,
$f:X\to\A^d_{K^+}:=\Spec\,K^+[T_1,\dots,T_d]$ an {\'e}tale
morphism. For every $r\in\N$ we consider the cartesian diagram
of schemes
$$
\xymatrix{
X_r \ar[r]^-{f_r} \ar[d]_{\psi_r} & \A^d_{K^+} \ar[d]^{\phi_r} \\
X \ar[r]^-f & \A^d_{K^+}
}$$
where $\phi_r$ is the morphism corresponding to the
$K^+$-algebra homomorphism
$$
\phi^\natural_r:K^+[T_1,\dots,T_d]\to K^+[T_1,\dots,T_d]
$$
defined by the rule: $T_j\mapsto T_j^{p^r}$ for $j=1,\dots,d$. For
every $r,s\in\N$ with $r\geq s$, $\psi_r$ factors through an obvious
$S$-morphism $\psi_{rs}:X_r\to X_s$, and the collection of the
schemes $X_r$ and transition morphisms $\psi_{sr}$ gives rise to an
inverse system $\underline X:=(X_r~|~r\in\N)$, whose inverse limit
is representable by an $S$-scheme $X_\infty$ (\cite[Ch.IV,
Prop.8.2.3]{EGAIV-3}). Let $g_\infty:X_\infty\to S$ (resp.
$g_r:X_r\to S$ for every $r\in\N$) be the structure morphism, $x\in
g_\infty^{-1}(s)$ any point, $x_r\in X_r$ the image of $x$ and
$R_r:=\cO_{\!X_r,x_r}$ for every $r\in\N$. Clearly the colimit
$R_\infty$ of the inductive system $(R_r~|~r\in\N)$ is naturally
isomorphic to $\cO_{\!X_\infty,x}$. Moreover, notice that the
restriction $g_{r+1}^{-1}(s)\to g_r^{-1}(s)$ is a radicial morphism
for every $r\in\N$. It follows easily that the transition maps
$R_r\to R_{r+1}$ are finite; furthermore, by inspection one sees
that $\phi_r^\natural$ is flat and finitely presented, so $R_{r+1}$
is a free $R_r$-module of rank $p^d$, for every $r\in\N$. Hence, the
present situation is a special case of example \ref{ex_flat-case},
and therefore conditions (a) and (b) hold if we choose the sequence
of integers $(d_i~|~i\in\N)$ with
$$
d_i:=p^{id}/[\kappa(x_i):\kappa(x_0)]
\qquad\text{for every $i\in\N$}.
$$
The foregoing discussion then yields a well-behaved notion
of normalized length for arbitrary $R_\infty$-modules supported
at $\{s(R_\infty)\}$.
\end{example}

\begin{example} In the situation of example \ref{ex_smooth},
it is easy to construct $R^a_\infty$-modules $M\neq 0$
such that $\lambda_\infty(M)=0$. For instance, for $d:=1$,
let $x\in X_\infty$ be the point of the special fibre
where $T_1=0$; then we may take $M:=R_\infty/I$, where $I$ is the
ideal generated by a non-zero element of $\fm_K$ and by the radical
of $T_1R_\infty$. The verification shall be left to the reader.
\end{example}

\subsection{Finite group actions on almost algebras}
\label{sec_fin-groups-quot}
In this section we fix a basic setup $(V,\fm)$ such that
$\tilde\fm:=\fm\otimes_V\fm$ is a flat $V$-module (see
\cite[\S2.1.1]{Ga-Ra}), and we consider some descent problems
for $V^a$-algebras endowed with a finite group of automorphisms.
Hence, the results below overlap with those of \cite[\S4.5]{Ga-Ra}.

\sset\subsubsection{}
Let $G$ be a finite group, $A$ a $V^a$-algebra, and let
$S:=\Spec\,V^a$, $X:=\Spec\,A$. A {\em right action of\/ $G$ on $X$\/}
is a group homomorphism :
$$
\rho:G\to\Aut_{V^a\Alg}(A)
$$
from $G$ to the group of automorphisms of $A$. Let $G_S$ be the
affine group $S$-scheme defined by $G$; hence every $g\in G$
determines a section $g_S:S\to G_S$ of the structure morphism
$G_S\to S$, and the resulting morphism:
$$
\coprod_{g\in G}g_S:S\amalg\cdots\amalg S\to G_S
$$
is an isomorphism of $S$-schemes. Then $\rho$ can be also regarded as
a right action of $G_S$ on $X$, as defined in \cite[\S3.3.6]{Ga-Ra}.
Especially, $\rho$ induces morphisms of $S$-schemes :
$$
\partial_i:X\times G:=X\times_SG_S\to X\qquad i=0,1
$$
as in {\em loc. cit.}, and we may define a {\em $G$-action on an
$A$-module $M$\/} (covering the given action of $G$ on $X$) as a
morphism of quasi-coherent $\cO_{\!X\times G}$-modules :
$$
\beta:\partial^*_0M\to\partial^*_1M
$$
fulfilling the conditions of \cite[\S3.3.7]{Ga-Ra}. One also says
that $(M,\beta)$ is a {\em $G$-equivariant $A$-module}. We denote by
$A[G]\Mod$ the category of all $G$-equivariant $A$-modules and
$G$-equivariant $A$-linear morphisms. Notice that $A[G]\Mod$ is
an abelian tensor category : indeed, for any two objects $(M,\beta)$,
$(M,\beta')$ we may set $(M,\beta)\otimes_A(M',\beta'):=
(M\otimes_AM',\beta\otimes_{\cO_{\!X\times G}}\beta')$.

\sset\subsubsection{}
Likewise, if $B$ is any $A$-algebra, a $G$-action on $B$ is a
morphism $\beta:\partial^*_0B\to\partial^*_1B$ of quasi-coherent
$\cO_{\!X\times G}$-algebras, such that the pair $(B,\beta)$
is a $G$-equivariant $A$-module. We say that $(B,\beta)$ is a
{\em $G$-equivariant $A$-algebra}, and we denote by $A[G]\Alg$
the category of such pairs, with $G$-equivariant morphisms
of $A$-algebras. One verifies easily that the datum $(B,\beta)$
is the same as a morphism $\psi:A\to B$ of $V^a$-algebras,
together with a $G$-action $\rho_B:G\to\Aut_{V^a\Alg}(B)$ on the
affine scheme $\Spec\,B$, such that the diagram :
$$
\xymatrix{ A \ar[r]^\psi \ar[d]_{\rho(g)} & B \ar[d]^{\rho_B(g)} \\
           A \ar[r]^\psi & B
}$$
commutes for every $g\in G$. We shall also consider the full subcategory
$A[G]\Alg_\mathrm{fl}$ (resp. $A[G]\wEt$, resp. $A[G]\Et$, resp.
$A[G]\Et_\mathrm{afp}$) of all such pairs, where $B$ is a flat
(resp. weakly {\'e}tale, resp. {\'e}tale, resp. {\'e}tale and almost
finitely presented) $A$-algebra.

\sset\subsubsection{}\label{subsec_triv-action}
The {\em trivial $G$-action\/} on $X$ is the map $\rho$ with
$\rho(g)=\one_A$ for every $g\in G$; this is the same as saying that
$\partial_0=\partial_1$. If $G$ acts trivially on $X$, a
$G$-equivariant $A$-module $(M,\beta)$ is the same as a
group homomorphism $\bar\beta:G\to\Aut_A(M)$ from $G$ to the group of
$A$-linear automorphisms of $M$. Namely, for every $g\in G$,
one lets :
\set\begin{equation}\label{eq_ties}
\bar\beta(g):=(\one_X\times_Sg_S)^*\beta
\end{equation}
and conversely, to a given map $\bar\beta$ there corresponds
a unique pair $(M,\beta)$ such that \eqref{eq_ties} holds.

Under this correspondence, the {\em trivial $G$-action\/}
$\bar\beta_0$ (such that $\bar\beta(g)=\one_M$ for every $g\in\ G$)
corresponds to the identity morphism
$\beta_0:\partial^*_0M\isom\partial^*_1M$.

More generally, let $(M,\beta)$ be a $G$-action on an $A$-module $M$,
covering the trivial $G$-action on $X$; for every $A_*$-valued
character $\chi:G\to A_*^\times$ of $G$, we let
$$
M_\chi:=\bigcap_{g\in G}\Ker(\bar\beta(g)-\chi(g)\cdot\one_M).
$$
The restriction of $\bar\beta$ defines a $G$-action on $M_\chi$,
such that the monomorphism $M_\chi\subset M$ is $G$-equivariant.
In the special case where $\chi$ is the trivial character, we have
$M_\chi=(M,\beta)^G$, the largest $G$-equivariant $A$-submodule of
$(M,\beta)$ on which $\beta$ restricts to the trivial $G$-action
({\em i.e.} the submodule fixed by $G$). When the notation is
not ambiguous, we shall often just write $M^G$ instead of $(M,\beta)^G$.

\sset\subsubsection{}
If $H\subset G$ is a subgroup, any $G$-action $\rho$ on $X$ induces
by restriction an $H$-action $\rho_{|H}$; then the morphisms
$\partial_i:X\times H\to X$ are just the restrictions of the
corresponding morphisms for $G$ (under the natural closed immersion
$X\times H\to X\times G$).
Similarly, a $G$-action $\beta$ on an $A$-module $M$ induces
by restriction an $H$-action $\beta_{|H}$ on the same module.

Let $Y:=\Spec\,B$ be any affine $S$-scheme. The $G$-action $\rho$
induces a $G$-action $Y\times_S\rho$ on $Y\times_SX$; namely,
$Y\times_S\rho(g):=\one_Y\times_S\rho(g)$ for every $g\in G$.
In terms of the group scheme $G_S$, this is the action given
by the morphisms :
$$
\partial_{Y,i}:=
\one_Y\times_S\partial_i:(Y\times_SX)\times G\to(Y\times_SX)
\qquad\text{$i=0,1$}.
$$
Let $\pi_X:Y\times_SX\to X$ be the natural morphism; then every
$G$-equivariant $A$-module $(M,\beta)$ induces a $G$-equivariant
$B\otimes_{V^a}A$-module $\pi^*_X(M,\beta):=(\pi^*_XM,\pi^*_X\beta)$,
whose action covers the $G$-action $Y\times_S\rho$ on $Y\times_SX$.

\sset\subsubsection{}\label{subsec_co-equalizers}
Let us set :
$$
\xymatrix{X/G:=\mathrm{Coequal}(X\times G\ar@<.5ex>[r]^-{\partial_1}
\ar@<-.5ex>[r]_-{\partial_1} & X) \quad\text{and}\quad
X^{\langle g\rangle}:=\Equal(X\ar@<.5ex>[r]^-{\rho(g)}
\ar@<-.5ex>[r]_-{\one_X} & X)}\quad\text{for every $g\in G$}.
$$
In other words, $X/G=\Spec\,A^G$, the subalgebra fixed by $G$, and
$X^{\langle g\rangle}$ is the closed subscheme fixed by the subgroup
${\langle g\rangle}\subset G$ generated by $g$. Thus,
$X^{\langle g\rangle}$ is the spectrum of a quotient $A/I_g$ of $A$,
where $I_g\subset A$ is the ideal generated by the almost elements
of the form $a-\rho(g)(a)$, for every $a\in A_*$.

Let $\pi:X\to X/G$ be the natural morphism, $N$ any quasi-coherent
$\cO_{\!X/G}$-module; the pull-back $\pi^*N$ is the quasi-coherent
$\cO_{\!X}$-module $A\otimes_{A^G}N$. By construction, there is a
natural isomorphism
$\beta_N:\partial_0^*(\pi^*N)\isom\partial_1^*(\pi^*N)$ (deduced
from the natural isomorphisms of functors
$\partial_i^*\circ\pi^*\simeq(\pi\circ\partial_i)^*$, for $i=0,1$),
and one verifies easily that $\beta_N$ is a $G$-action on $\pi^*N$.

Moreover, let $i_g:X^{{\langle g\rangle}}\to X$ be the natural
closed immersion; if $M$ is any $A$-module, then $i^*_gM$ is the
$A/I_g$-module $M/I_gM$. Especially, take $M:=\pi^*N$; we notice that
the restriction $\rho_{|\langle g\rangle}$ of the given action
$\rho$, induces the trivial $\langle g\rangle$-action on
$X^{\langle g\rangle}$, and directly from the construction,
we see that the natural action $\beta_{N|\langle g\rangle}$
restricts to the trivial $\langle g\rangle$-action on $i_g^*M$.

We are thus led to the :
\begin{definition}\label{def_horizontal}
Let $G$ be a finite group, $A$ a $V^a$-algebra, $\rho$ a $G$-action
on $X:=\Spec\,A$.
\begin{enumerate}
\item
The category $A[G]\Mod_\hor$ of
{\em $A$-modules with horizontal $G$-action\/} is the full
subcategory of $A[G]\Mod$ consisting of all pairs $(M,\beta)$
subject to the following condition. For every $g\in G$, the
restriction $\beta_{|\langle g\rangle}$ induces the trivial
action on $i^*_gM$.
\item
We denote by $A[G]\Mod_\horf$ the full subcategory of $A[G]\Mod_\hor$
consisting of the pairs $(M,\beta)$ as above, such that $M$ is a flat
$A$-module.
\item
Likewise, we denote by $A[G]\Alg_\hor$ (resp. $A[G]\Alg_\horf$, resp.
$A[G]\wEt_{\hor}$, resp. $A[G]\Et_{\hor}$,
resp. $A[G]\Et_{\hor.\mathrm{afp}}$) the full subcategory of
$A[G]\Alg$ (resp. $A[G]\Alg_\mathrm{fl}$, resp. $A[G]\wEt$,
resp. $A[G]\Et$, resp. $A[G]\Et_\mathrm{afp}$) consisting of
all pairs $(B,\beta)$ which are horizontal, when regarded as
$G$-equivariant $A$-modules.
\end{enumerate}
\end{definition}

Notice that the tensor product of horizontal modules (resp. algebras)
is again horizontal. By the foregoing, the rule
$N\mapsto(\pi^*N,\beta_N)$ defines a functor :
\set\begin{equation}\label{eq_horizontal}
A^G\Mod\to A[G]\Mod_\hor.
\end{equation}
On the other hand, if $(M,\beta)$ is any $G$-equivariant $A$-module,
the pair $\pi_*(M,\beta):=(\pi_*M,\pi_*\beta)$ may be regarded as
a $G$-action on $\pi_*M$, covering the trivial $G$-action on $\Spec\,A^G$,
{\em i.e.} a group homomorphism $G\to\Aut_{A^G}(M)$ (see
\eqref{subsec_triv-action}). One verifies easily that the functor
$N\mapsto(\pi^*N,\beta_N)$ is left adjoint to the functor
$A[G]\Mod\to A^G\Mod$ : $(M,\beta)\mapsto\pi_*(M,\beta)^G$ (details left to
the reader). Hence, also \eqref{eq_horizontal} admits a right adjoint,
given by the same rule. Similar assertions hold for the analogous
functors :
\set\begin{equation}\label{eq_same-for-algs}
A^G\Alg\to A[G]\Alg_\hor
\end{equation}
and the variants considered in definition \ref{def_horizontal}(iii).

\begin{lemma}\label{lem_horizontal} In the situation of definition
 {\em\ref{def_horizontal}}, suppose furthermore that the order $o(G)$
of $G$ is invertible in $A_*$, and let $(M,\beta)$ be any
$G$-equivariant $A$-module. Then :
\begin{enumerate}
\item
For every $A^G_*$-valued character $\chi:G\to(A^G_*)^\times$, the natural
$A^G$-linear monomorphism $\pi_*M_\chi\to\pi_*M$ admits a $G$-equivariant
$A^G$-linear right inverse $\pi_*M\to\pi_*M_\chi$.
\item
For every $A^G$-module $N$, the unit of adjunction :
$$
\eps_N:N\to(A\otimes_{A^G}N)^G
$$
is an isomorphism.
\item
Let $Y:=\Spec\,B$ be any affine $S$-scheme; denote by
$\pi_X:Y\times_SX\to X$ and $\pi_{X/G}:Y\times_S(X/G)\to X/G$
the natural projections. Then the natural morphism :
$$
\pi^*_{X/G}(M,\beta)^G\to(\pi^*_X M,\pi^*_X\beta)^G
$$
is an isomorphism of $B\otimes_{V^a}A^G$-modules.
\end{enumerate}
\end{lemma}
\begin{proof} This is standard : for every $\chi$ as in (i), the
group algebra $A_*[G]$ admits the central idempotent :
\set\begin{equation}\label{eq_central-idemp}
e_\chi:=\frac{1}{o(G)}\cdot\sum_{g\in G}\chi(g)\cdot g
\end{equation}
and $M_\chi=e_\chi M$ for any $G$-equivariant $A$-module $M$. Especially,
we may take $M:=N\otimes_{A^G}A$, and $e_0$ the central idempotent
associated with the trivial character, in which case $e_0A=A^G$, and
$M=N\oplus(N\otimes_{A^G}(1-e_0)A)$, so all the claims follow easily.
\end{proof}

\begin{definition}\label{def_horiz-graded}
Let $\Gamma$ be any finite abelian group with neutral element
$0\in\Gamma$.
\begin{enumerate}
\item
A {\em $\Gamma$-graded $V^a$-algebra\/} is a pair
$\underline A:=(A,\gr_\bullet A)$ consisting of a $V^a$-algebra
$A$ and a decomposition $A=\bigoplus_{\chi\in\Gamma}\gr_\chi A$
as a direct sum of $V^a$-modules, such that :
$$
1\in\gr_0A_*\quad\text{and}\quad
\gr_\chi A\cdot\gr_{\chi'}A\subset\gr_{\chi+\chi'}A
\qquad\text{for every $\chi,\chi'\in\Gamma$}
$$
(where as usual $\gr_\chi A\cdot\gr_{\chi'}A$ denotes the
image of the restriction $\gr_\chi A\otimes_{V^a}\gr_{\chi'}A\to A$
of the multiplication morphism $\mu_A$).
Especially, $\gr_0A$ is a $V^a$-subalgebra of $A$, and  every
submodule $\gr_\chi A$ is a $\gr_0A$-module.
\item
A {\em $\Gamma$-graded $\underline A$-module\/} is a pair
$\underline N:=(N,\gr_\bullet N)$ consisting of an $A$-module
$N$ and a decomposition $N=\bigoplus_{\chi\in\Gamma}\gr_\chi N$
as a direct sum of $V^a$-modules, such that :
$$
\gr_\chi A\cdot\gr_{\chi'} N\subset\gr_{\chi+\chi'}N
\qquad\text{for every $\chi,\chi'\in\Gamma$}.
$$
Of course, a morphism of $\Gamma$-graded $\underline A$-modules
$\underline N\to\underline N':=(N',\gr_\bullet N')$ is an $A$-linear
morphism $N\to N'$ that respects the gradings.
\item
For every subgroup $\Delta\subset\Gamma$, let $J_\Delta\subset A$
be the graded ideal generated by
$\bigoplus_{\chi\notin\Delta}\gr_\chi A$. We say that $\underline N$
is {\em horizontal\/} if
$\gr_\chi(N/J_\Delta N):=\gr_\chi N/(\gr_\chi N\cap J_\Delta N)=0$
for every subgroup $\Delta\subset\Gamma$ and every $\chi\notin\Delta$.
\item
If $\Delta\subset\Gamma$ is any subgroup, $p:\Gamma\to\Gamma/\Delta$
the natural projection, and $\rho\in\Gamma/\Delta$ any element, we let :
$$
\gr^{\Gamma/\Delta}_\rho N:=\bigoplus_{\chi\in p^{-1}(\rho)}\gr_\chi N
$$
and set $\underline N_{\Gamma/\Delta}:=(N,\gr^{\Gamma/\Delta}_\bullet N)$.
Then $\underline A_{\Gamma/\Delta}$ is a $\Gamma/\Delta$-graded $V^a$-algebra,
and $\underline N_{\Gamma/\Delta}$ is a $\Gamma/\Delta$-graded
$\underline A_\Delta$-module. Moreover, let also
$\underline N_{|\Delta}$ be the pair consisting of
$N_{|\Delta}:=\gr^{\Gamma/\Delta}_0N$
together with its decomposition
$N_{|\Delta}=\bigoplus_{\chi\in\Delta}\gr_\chi N$; then
$\underline A_{|\Delta}$ is a $\Delta$-graded $V^a$-algebra, and
$\underline N_{|\Delta}$ is a $\Delta$-graded
$\underline A_{|\Delta}$-module.
\end{enumerate}
\end{definition}

\begin{proposition}\label{prop_several-years-on}
In the situation of definition {\em\ref{def_horiz-graded}},
let $\underline N$ be any $\Gamma$-graded $\underline A$-module.
Then the following conditions are equivalent :
\begin{enumerate}
\alphaenu
\item
$\underline N$ is horizontal.
\item
The natural morphism $A\otimes_{\gr_0A}\gr_0N\to N$ is an epimorphism.
\end{enumerate}
\end{proposition}
\begin{proof}(b)$\Rightarrow$(a): The assertion is obvious
for the $\Gamma$-graded $\underline A$-module consisting of
$A\otimes_{\gr_0A}\gr_0N$ and its natural grading deduced from
$\gr_\bullet A$; however any (graded) quotient of a horizontal
module is horizontal, hence the assertion follows also for
$\underline N$.

(a)$\Rightarrow$(b): We argue by induction on $o(\Gamma)$.
The first case is covered by the following :

\begin{claim}\label{cl_OK-for-prime-ord}
The proposition holds if $o(\Gamma)$ is a prime number.
\end{claim}
\begin{pfclaim} Indeed, suppose that $\underline N$ is horizontal.
We may replace $N$ by $N/(A\cdot\gr_0N)$, which is a horizontal
$\Gamma$-graded $\underline A$-module, when endowed with the
grading induced from $\gr_\bullet N$. Then $\gr_0N=0$, and we have
to show that $N=0$. By assumption, $N\subset J_{\{0\}}N$, {\em i.e.} :
$$
\gr_\chi N\subset
\sum_{\sigma\neq 0,\chi}\gr_\sigma A\cdot\gr_{\chi-\sigma}N
\qquad\text{for every $\chi\neq 0$}.
$$
For every $n>0$, the symmetric group $S_n$ acts on the set
$(\Gamma\setminus\{0\})^n$ by permutations; we let
$Q_n:=(\Gamma\setminus\{0\})^n/S_n$, the set of equivalence classes
under this action. For every
$\underline\sigma:=(\sigma_1,\dots,\sigma_n)\in Q_n$,
let $|\underline\sigma|:=\sum_{i=1}^n\sigma_i$. By an easy induction,
it follows that :
$$
\gr_\chi N\subset\sum_{\underline\sigma\in Q_n}
\gr_{\sigma_1}A\cdots\gr_{\sigma_n}A\cdot
\gr_{\chi-|\underline\sigma|}N.
$$
Since every element of $\Gamma\setminus\{0\}$ generates $\Gamma$,
it is also clear that there exists $n\in\N$ large enough such
that every sequence $\underline\sigma$ in $Q_n$ admits a subsequence,
say $\underline\tau:=(\tau_1,\dots,\tau_m)$ for some $m\leq n$, with
$|\underline\tau|=\chi$ (details left to the reader). Up to a
permutation, we may assume that $\underline\tau$ is the final
segment of $\underline\sigma$; then we have :
$$
\gr_{\sigma_1}A\cdots\gr_{\sigma_n}A\cdot
\gr_{\chi-|\underline\sigma|}N\subset
\gr_{\sigma_1}A\cdots\gr_{\sigma_{n-m}}A\cdot
\gr_0N=0
$$
whence the claim.
\end{pfclaim}

Next, suppose that the assertion is already known for every subgroup
$\Gamma'\subset\Gamma$, every graded $\Gamma'$-algebra $\underline B$,
and every $\Gamma'$-graded $\underline B$-module $\underline P$.
We choose a subgroup $\Gamma'\subset\Gamma$ such that
$(\Gamma:\Gamma')$ is a prime number. We shall use the following :

\begin{claim}\label{cl_Russell}
Let $G$ be a group, $0\in G$ the neutral element, $H\subset G$
a subgroup. For every subgroup $L\subset G$ with $H\cap L=\{0\}$,
choose an element $g_L\in G\setminus L$; denote by $S$ the subgroup
generated by all these elements $g_L$. Then $S\cap H\neq\{0\}$.
\end{claim}
\begin{pfclaim} Indeed, if $S\cap H=\{0\}$, we would have $g_S\in S$,
a contradiction.
\end{pfclaim}

\begin{claim}\label{cl_horiz-step}
Suppose that $\underline N$ is horizontal. Then
$\underline N_{|\Gamma'}$ is a horizontal $\Gamma'$-graded
$\underline A_{|\Gamma'}$-module.
\end{claim}
\begin{pfclaim} For any given subgroup $\Delta\subset\Gamma'$,
let $J'_\Delta\subset A_{|\Gamma'}$ be the ideal generated
by $\bigoplus_{\chi\in\Gamma'\setminus\Delta}\gr_\chi A$; have
to show that $\gr_\chi(N/J'_\Delta N)=0$ for every
$\chi\in\Gamma'\setminus\Delta$. To this aim, we may replace
$\Gamma$, $\Gamma'$, $\underline A$, $\underline A_{|\Gamma'}$,
$\underline N$ and $\underline N_{|\Gamma'}$,
by respectively $\Gamma/\Delta$, $\Gamma'/\Delta$,
$\underline A_{\Gamma/\Delta}$,
$(\underline A_{|\Gamma'})_{\Gamma'/\Delta}$,
$\underline N_{\Gamma/\Delta}$ and
$(\underline N_{|\Gamma'})_{\Gamma'/\Delta}$, which allows
to assume that $\Delta=\{0\}$. In this case, we have to show
that :
\set\begin{equation}\label{eq_reduce-to--zero}
N_{|\Gamma'}=\gr_0N+J'_{\{0\}}\cdot N_{|\Gamma'}.
\end{equation}
Notice that $A_{|\Gamma'}\cdot \gr_0N=\gr_0N+J'_{\{0\}}\cdot \gr_0N$; the
quotient $P:=N/(A\cdot\gr_0N)$ carries a unique $\Gamma$-grading
$\gr_\bullet P$ such that the projection
$\underline N\to\underline P:=(P,\gr_\bullet P)$ is a morphism
of $\Gamma$-graded $\underline A$-modules (namely,
$\gr_\chi P:=N_\chi/\gr_\chi A\cdot\gr_0N$ for every $\chi\in\Gamma$).
Moreover, $P$ is horizontal, $\gr_0P=0$, and
\eqref{eq_reduce-to--zero} is equivalent to :
$P_{|\Gamma'}=J'_{\{0\}}\cdot P_{|\Gamma'}$.
Therefore we may replace $\underline N$ by $\underline P$, and
assume from start that $\gr_0N=0$. Similarly, notice that
$J'_{\{0\}}\cdot N$ is a $\Gamma$-graded $\underline A$-submodule
of $\underline N$, and the pair $\underline Q$ consisting of
$Q:=N/(J'_{\{0\}}\cdot N)$ and its quotient grading, is horizontal.
Moreover, $\gr_\chi Q=\gr_\chi(N_{|\Gamma'}/J'_{\{0\}}N_{|\Gamma'})$,
for every $\chi\in\Gamma'$. Hence, we may replace $N$ by $Q$,
which allows to assume as well that
\set\begin{equation}\label{eq_assume-two}
J'_{\{0\}}N_{|\Gamma'}=0
\end{equation}
in which case, we are reduced to showing that $N=0$. Now, by
assumption, for every subgroup $H\subset\Gamma$ we have :
$$
N=N_{|H}+J_HN.
$$
Let $\cC$ be the set of all subgroups $H\subset\Gamma$ with
$H\cap\Gamma'=\{0\}$; we deduce that :
\set\begin{equation}\label{eq_omit-a-few}
\gr_\chi N\subset
\sum_{\sigma\notin H}\gr_\sigma A\cdot\gr_{\chi-\sigma}N
\qquad\text{for every $\chi\in\Gamma'\setminus\{0\}$ and every $H\in\cC$}.
\end{equation}
Moreover, in view of \eqref{eq_assume-two} we may omit from
the sum in \eqref{eq_omit-a-few} all the elements $\sigma$ that lie in
$\Gamma'$; for the remaining elements we have
$\chi-\sigma\notin\Gamma'$, hence we may apply claim
\ref{cl_OK-for-prime-ord} to the horizontal $\Gamma/\Gamma'$-graded
$\underline A_{\Gamma/\Gamma'}$-module $\underline N_{\Gamma/\Gamma'}$,
to deduce that :
$$
\gr_{\chi-\sigma}N\subset\sum_{\delta\in\Gamma'\setminus\{0\}}
\gr_{\chi-\sigma-\delta}A\cdot\gr_\delta N.
$$
Therefore :
$$
\gr_\chi N\subset\sum_{\sigma\notin H\cup\Gamma'}
\sum_{\delta\in\Gamma'\setminus\{0\}}\gr_\sigma A\cdot
\gr_{\chi-\sigma-\delta}A\cdot\gr_\delta N
\qquad\text{for every $\chi\in\Gamma'\setminus\{0\}$ and every $H\in\cC$}.
$$
However -- again due to \eqref{eq_assume-two} -- we may omit from
this sum all the terms corresponding to the pairs $(\sigma,\delta)$
with $\chi\neq\delta$, hence we conclude that :
$$
\gr_\chi N\subset\sum_{\sigma\notin H\cup\Gamma'}
\gr_\sigma A\cdot\gr_{-\sigma}A\cdot\gr_\chi N
\qquad\text{for every $\chi\in\Gamma'\setminus\{0\}$ and every
$H\in\cC$}.
$$
Denote by $\Sigma$ the set of all mappings
$\underline\sigma:\cC\to\Gamma$ such that
$\underline\sigma(H)\notin H\cup\Gamma'$ for every $H\in\cC$.
Moreover, for every $\underline\sigma\in\Sigma$, set :
$$
B_{\underline\sigma}:=\prod_{H\in\cC}
\gr_{\underline\sigma(H)}A\cdot\gr_{-\underline\sigma(H)}A
$$
(this is an ideal of $\gr_0A$) and let
$S_{\underline\sigma}\subset\Gamma$ be the subgroup generated
by the image of $\underline\sigma$. By an easy induction we deduce
that :
$$
\gr_\chi N\subset\sum_{\underline\sigma\in\Sigma}
B^n_{\underline\sigma}\cdot\gr_\chi N
\qquad\text{for every $n>0$}.
$$
By claim \ref{cl_Russell}, for every $\underline\sigma\in\Sigma$
we may find
$\gamma(\underline\sigma)\in S_{\underline\sigma}\cap\Gamma'\setminus\{0\}$.
On the other hand, since $\Gamma$ is finite and abelian, it is
easy to verify that there exists $n\in\N$ large enough such that
$$
B^n_{\underline\sigma}\subset
\gr_{\gamma(\underline\sigma)}A\cdot\gr_{-\gamma(\underline\sigma)}A
\qquad\text{for every $\underline\sigma\in\Sigma$}
$$
(details left to the reader). But \eqref{eq_assume-two} implies that
$\gr_{\gamma(\underline\sigma)}A\cdot\gr_\chi N=0$ whenever
$\chi\in\Gamma'\setminus\{0\}$, so the claim follows.
\end{pfclaim}

To conclude, we apply first claim \ref{cl_OK-for-prime-ord} to the
horizontal $\Gamma/\Gamma'$-graded $\underline A_{\Gamma/\Gamma'}$-module
$\underline N_{\Gamma/\Gamma'}$, to see that $N_{|\Gamma/\Gamma'}$
generates the $A$-module $N$, and then claim \ref{cl_horiz-step}
-- together with our inductive assumption -- to deduce that $\gr_0N$
generates that $A_{|\Gamma/\Gamma'}$-module $N_{|\Gamma/\Gamma'}$.
The proposition follows.
\end{proof}

\sset\subsubsection{}\label{subsec_pure}
Recall that a morphism $M\to N$ of $A$-modules is said to be
{\em pure\/} if the natural morphism $Q\otimes_AM\to Q\otimes_AN$
is a monomorphism for every $A$-module $Q$. A morphism $A\to B$
of $V^a$-algebras is called {\em pure\/} if it is pure when
regarded as a morphism of $A$-modules.

\begin{lemma}\label{lem_pure-desc-flat}
Let $f:A\to B$ be a pure morphism of $V^a$-algebras, $M$ an
$A$-module. Then :
\begin{enumerate}
\item
If $B\otimes_AM$ is a flat $B$-module, then $M$ is a flat $A$-module
({\em i.e.} $f$ descends flatness).
\item
If $B\otimes_AM$ is almost finitely generated (resp. almost finitely
presented) as a $B$-module, then $M$ is almost finitely generated
(resp. almost finitely presented) as an $A$-module.
\end{enumerate}
\end{lemma}
\begin{proof}(i): To begin with, we remark :

\begin{claim}\label{cl_about-pure}
Let $R$ be any $V$-algebra such that $A=R^a$.
\begin{enumerate}
\item
A morphism $\phi:M_1\to M_2$ of $A$-modules is pure if and only if
the same holds for the induced morphism $\phi_!:M_{1!}\to M_{2!}$ of
$R$-modules.
\item
A morphism $g:A\to B'$ of $V^a$-algebras is pure if and only
if the same holds for the induced morphism $g_{!!}:A_{!!}\to B'_{!!}$.
\end{enumerate}
\end{claim}
\begin{pfclaim} (i): Suppose that $\phi$ is pure, and let $Q$
be any $R$-module. Then $Q\otimes_RM_{i!}\simeq(Q^a\otimes_AM_i)_!$
for $i=1,2$. Since the functor $M\mapsto M_!$ is exact
(\cite[Cor.2.2.24(i)]{Ga-Ra}), we deduce that $\phi_!$ is pure.
The converse is easy, and shall be left to the reader.

(ii): From (i) we already see that $g$ is pure whenever $g_{!!}$
is. Next, suppose that $g$ is pure; we may assume that $R=A_{!!}$,
and then (i) says that $g_!$ is a pure morphism of
$A_{!!}$-modules. However, the natural diagram of $A_{!!}$-modules :
$$
\xymatrix{ A_! \ar[r] \ar[d] & B'_! \ar[d] \\
           A_{!!} \ar[r] & B'_{!!}
}$$
is cofibred; since tensor products are right exact functors,
the claim follows.
\end{pfclaim}

The assertion now follows from claim \ref{cl_about-pure}(ii) and
\cite[Partie II, Lemme 1.2.1]{Gr-Ra}.

(ii): Suppose first that $B\otimes_AM$ is an almost finitely
generated $B$-module. The assumption on $f$ implies that
$\Ann_A(B\otimes_AM)\subset\Ann_A M$; then
\cite[Rem.3.2.26(i)]{Ga-Ra} shows that $M$ is an almost finitely
generated $A$-module.

Finally, we suppose that $B\otimes_AM$ is an almost finitely
presented $B$-module, and we wish to show that $M$ is an almost
finitely presented $A$-module. To this aim, let $\phi:N\to N'$
be a morphism of $A$-modules. The assumption on $f$ implies
that the natural morphism
$\Ker\,\phi\to\Ker(\one_B\otimes_A\phi)$ is a monomorphism;
especially :
$$
\Ann_A(\Ker(\one_B\otimes_A\phi))\subset\Ann_A(\Ker\,\phi).
$$
Then one may easily adapt the proof of \cite[Lemma 3.2.25(iii)]{Ga-Ra},
to derive the assertion.
\end{proof}

\begin{theorem}\label{th_horiz}
In the situation of definition {\em\ref{def_horizontal}}, suppose
furthermore that $G$ is abelian and the order $o(G)$ of $G$ is
invertible in $A_*$. Then the following holds :
\begin{enumerate}
\item
Let $M$ be any $G$-equivariant $A$-module. The $G$-action
on $M$ is horizontal if and only if the counit of adjunction :
$$
\eta_M:A\otimes_{A^G}M^G\to M
$$
is an epimorphism (of $A$-modules).
\item
The functor \eqref{eq_horizontal} restricts to an equivalence
on the full subcategory of flat $A^G$-modules :
\set\begin{equation}\label{eq_flat-horizon}
A^G\Mod_\fl\isom A[G]\Mod_\horf.
\end{equation}
\end{enumerate}
\end{theorem}
\begin{proof}(i): For $m:=o(G)$, let $\bmu_m\subset\bar\Q{}^\times$
be the group of $m$-th roots of $1$, and set
$B:=V^a[\bmu_m]:=(V[T]/(T^m-1))^a$. Since $B$ is a faithfully flat
$V^a$-algebra, lemma \ref{lem_horizontal}(iii) allows to replace
$A$ by $B\otimes_{V^a}A$ and $M$ by $B\otimes_{V^a}M$, and therefore
we may assume from start that $\bmu_m\subset(A^G_*)^\times$.
Set $\Gamma:=\Hom_\Z(G,\bmu_m)$. For every $\chi\in\Gamma$, let
$e_\chi\in A_*[G]$ be the central idempotent defined as in
\eqref{eq_central-idemp}. A standard calculation shows that :
$$
\sum_{\chi\in\Gamma}e_\chi=1.
$$
Hence, every $G$-equivariant $A^G$-module $(N,\beta)$ admits
the $G$-equivariant decomposition :
\set\begin{equation}\label{eq_dec-by-idemp}
N\simeq\bigoplus_{\chi\in\Gamma}N_\chi.
\end{equation}
Especially, $A=\bigoplus_{\chi\in\Gamma}A_\chi$, and clearly the
datum $\underline A$ consisting of $A$ and its decomposition, is
a $\Gamma$-graded $V^a$-algebra.
Furthermore, the datum $\underline N$ consisting of $N$ and its decomposition
\eqref{eq_dec-by-idemp} is a $\Gamma$-graded $\underline A$-module.

\begin{claim}
The functor $(N,\beta)\mapsto\underline N$ is an equivalence from
$A[G]\Mod$ to the category of $\Gamma$-graded $\underline A$-modules.
\end{claim}
\begin{pfclaim} Indeed, if $Q$ is a $\Gamma$-graded
$\underline A$-module, we may define a $G$-action on $Q$ by requiring
that $Q_\chi=\gr_\chi Q$, the $\chi$-graded direct summand of $Q$.
This gives a quasi-inverse for the functor $(N,\beta)\mapsto\underline N$.
(Details left to the reader.)
\end{pfclaim}

\begin{claim}\label{cl_horiz-both}
Let $(N,\beta)$ be any $G$-equivariant $A$-module, and $\underline N$
its associated $\Gamma$-graded $\underline A$-module.
The following conditions are equivalent :
\begin{enumerate}
\alphaenu
\item
$(N,\beta)$ is horizontal.
\item
$\underline N$ is horizontal.
\end{enumerate}
\end{claim}
\begin{pfclaim} For every $g\in G$, let
$\Delta(g)\subset\Gamma$ be the subgroup consisting of all
$\chi\in\Gamma$ such that $\chi(g)=1$; a direct inspection of the
definitions shows that $J_{\Delta(g)}$ is the ideal $I_g$, as defined
in \eqref{subsec_co-equalizers}, and $(N,\beta)$ is horizontal if and
only if $N_\chi/(N_\chi\cap I_gN)=0$ for every $g\in G$ and every
$\chi\notin\Delta(g)$.
This already shows that (b)$\Rightarrow$(a); it also shows
that condition (b) holds for the subgroups $\Delta(g)$, when
$(N,\beta)$ is horizontal.
However, every subgroup of $\Gamma$ can be written in the form
$\Delta=\Delta(g_1)\cap\cdots\cap\Delta(g_n)$, for appropriate
$g_1,\dots,g_n\in G$, and then
$J_{\Delta(g_1)}+\cdots+J_{\Delta(g_n)}\subset J_\Delta$, hence
(b) follows for all subgroups.
\end{pfclaim}

Assertion (i) now follows from claim \ref{cl_horiz-both} and
proposition \ref{prop_several-years-on}.

(ii): Let $(M,\beta)$ be any object of $A[G]\Mod_\horf$,
and denote by $L$ the kernel of the counit $\eta_M$. Since $M$ is
a flat $A$-module and $\eta_M$ is an epimorphism by (i), it follows
that $i_g^*L$ is the kernel of $i^*_g\eta_M$, for every
$g\in G$ (notation of \eqref{subsec_co-equalizers}). Since the
category $A[G]\Mod_\hor$ is abelian, we deduce that the natural
$G$-action on $L$ is horizontal, and then (i) says that $L$ is
generated by $L^G$. But lemma \ref{lem_horizontal}(ii) easily implies
that $L^G=0$, so $\eta_M$ is an isomorphism. Next, letting $M:=A$
in lemma \ref{lem_horizontal}(i), we deduce easily that the natural
morphism $A^G\to A$ is pure, hence $M^G$ is a flat $A^G$-module,
by lemma \ref{lem_pure-desc-flat}(i). Now the assertion follows from
lemma \ref{lem_horizontal}(ii) and \cite[Prop.3.4.3]{Bor}.
\end{proof}

\begin{corollary}\label{cor_horiz}
In the situation of theorem {\em\ref{th_horiz}}, the following holds :
\begin{enumerate}
\item
The functor \eqref{eq_flat-horizon} restricts to an equivalence
from the subcategory of flat, almost finitely generated (resp. almost
finitely presented) $A^G$-modules, onto the subcategory of
$G$-equivariant, flat, horizontal and almost finitely generated
(resp. almost finitely presented) $A[G]$-modules.
\item
The functor \eqref{eq_same-for-algs} restricts to an equivalence :
$$
A^G\Alg_\mathrm{fl}\to A[G]\Alg_\horf
$$
and likewise for the subcategories of weakly {\'e}tale (resp. {\'e}tale,
resp. {\'e}tale and almost finitely presented) algebras.
\end{enumerate}
\end{corollary}
\begin{proof} (i) follows from theorem \ref{th_horiz}(ii),
lemma \ref{lem_pure-desc-flat}(ii), and the fact that the morphism
$A^G\to A$ is pure.

(ii): The assertion concerning $A^G\Alg_\mathrm{fl}$
is an immediate consequence of theorem \ref{th_horiz}. Next, let
$(B,\beta)$ be an object of $A[G]\wEt$; by the foregoing, $B$
descends to a flat $A^G$-algebra $B^G$ with a $G$-equivariant
isomorphism : $B\simeq A\otimes_{A^G}B^G$. However, on the one
hand, $B$ is -- by assumption -- a flat $B\otimes_AB$-algebra, and
on the other hand, $B\otimes_AB$ underlies a flat, horizontal $A[G]$-algebra
with $(B\otimes_AB)^G\simeq B^G\otimes_{A^G}B^G$; theorem
\ref{th_horiz} then says that $B^G$ is a flat
$B^G\otimes_{A^G}B^G$-algebra, whence the assertion for $A^G\wEt$.
Next, since an almost finitely generated module is almost projective
if and only if it is flat and almost finitely presented
(\cite[Prop.2.4.18]{Ga-Ra}), the assertion for $A^G\Et$ follows from
the same assertion for $A^G\wEt$ and lemma \ref{lem_pure-desc-flat}(ii).
Likewise, the assertion for {\'e}tale almost finitely presented
$A^G$-algebras follows from the assertion for $A^G\Et$ and lemma
\ref{lem_pure-desc-flat}(ii).
\end{proof}

\begin{remark} In case the $G$-action on $X$ is {\em free\/},
{\em i.e.} when $(\partial_0,\partial_1):X\times G\to X\times_SX$
is a monomorphism, corollary \ref{cor_horiz} also follows from
\cite[Prop.4.5.25]{Ga-Ra}.
\end{remark}

\subsection{Almost Witt vectors}
In this section we show that the construction of the ring
of Witt vectors descends to almost rings, and we study some
properties of the resulting functor of {\em almost
Witt vectors}. We begin with some general observation
concerning liftings of basic setups (in the sense of
\cite[\S2.1.1]{Ga-Ra}) along ring homomorphisms. These
preliminaries shall then be applied to find natural
basic setups on truncated rings of Witt vectors of
arbitrary rings $R$, by lifting given setups on $R$
along the $0$-th ghost map.

\begin{lemma}\label{lem_lift-a-struct-over-nilp}
Let $A$ be a ring, $I\subset A$ a nilpotent ideal,
$\pi:A\to A_0:=A/I$ the projection, and $\fm_0\subset A_0$
an ideal such that $\fm_0=\fm_0^2$. We have :
\begin{enumerate}
\item
There exists a unique ideal $\fm\subset A$ such that
$\fm^2=\fm$ and $\pi(\fm)=\fm_0$.
\item
$\fm$ is the smallest of the ideals $J\subset A$ such that
$\pi(J)=\fm_0$.
\item
$\fm$ fulfills condition  $(\bB)$ of\/ \cite[\S2.1.6]{Ga-Ra}
if and only if the same holds for $\fm_0$.
\end{enumerate}
\end{lemma}
\begin{proof}(i): Say that $I^N=0$ for some $N\in\N$, and
set $J:=\pi^{-1}\fm_0$. We remark :

\begin{claim}\label{cl_here-is-one}
$J^{2N+1}=J^{2N+2}$.
\end{claim}
\begin{pfclaim} Indeed, let $a_1,\dots,a_{2N+1}\in J^{2N+1}$;
we may find an integer $r\geq 0$ and elements
$b_{1i},c_{1i},\dots,b_{ri},c_{ri}\in\fm_0$ for every
$i=1,\dots,2N+1$, such that $\pi(a_i)=\sum_{j=1}^rb_{ji}c_{ji}$.
Lift each $b_{ji}$ and each $c_{ji}$ to elements
$b'_{ji},c'_{ji}\in J$; it follows that
$x_i:=a_i-\sum_{j=1}b'_{ji}c'_{ji}\in I$ for every
$i=1,\dots,2N+1$. Thus :
$$
\prod_{i=1}^{2N+1}a_i=\prod_{i=1}^{2N+1}
\Bigl(x_i+\sum_{j=1}^rb'_{ji}c'_{ji}\Bigr)
\in J^{2(N+1)}+I^N=J^{2N+2}
$$
whence the claim.
\end{pfclaim}

Set $\fm:=J^{2N+1}$; claim \ref{cl_here-is-one} implies that
$\fm=\fm^2$, and clearly $\pi(\fm)=\fm_0^{2N+1}=\fm$. Next,
let $\fm'\subset A$ be any ideal such that $\fm'^2=\fm'$
and $\pi(\fm')=\fm_0$; it follows that
$$
\fm'=\fm'^{2N+1}\subset J^{2N+1}=(\fm'+I)^{2N+1}
\subset\fm'^{N+1}=\fm'
$$
whence $\fm'=\fm$.

(ii): Let $J\subset A$ be any ideal such that $\pi(J)=\fm_0$;
it follows that
$$
\fm=\fm^N\subset(J+I)^N\subset J+I^N=J.
$$

(iii): It is clear that if $\fm$ fulfills condition
$(\bB)$, then the same holds for $\fm_0$. Conversely,
if condition  $(\bB)$ holds for $\fm_0$, then for every
integer $k>1$ the system of elements $(x^k~|~x\in\fm)$
generates an ideal $J\subset A$ such that $\pi(J)=\fm_0$
and $J\subset\fm$. By (ii), we must then have $J=\fm$,
which shows that condition $(\bB)$ holds for $\fm$.
\end{proof}

Let us also point out the following result,
which shall not be needed in the sequel :

\begin{lemma} Let $(A,\fm)$ be a basic setup and $f:B\to A$
a ring homomorphism such that :
\begin{enumerate}
\alphaenu
\item
$p^kA=0$ for some integer $k\in\N$.
\item
$\fm$ fulfills condition $(\bB)$ of\/ \cite[\S2.1.6]{Ga-Ra}.
\item
$f\otimes_\Z\F_p$ is invertible up to $\Phi^n$ for some
$n\in\N$, in the sense of\/ \cite[Def.3.5.8]{Ga-Ra}.
\end{enumerate}
Then there exists a unique basic setup $(B,\fn)$ with
$\fn A=\fm$ and with $\fn$ fulfilling condition $(\bB)$.
\end{lemma}
\begin{proof} Using lemma \ref{lem_lift-a-struct-over-nilp},
we are easily reduced to the case where $A$ and $B$ are
$\F_p$-algebras, and $f$ is invertible up to $\Phi^n$.
The latter means that there exists a morphism $g:A\to B$
such that $g\circ f=\Phi^n_B$ and $f\circ g=\Phi^n_A$.

Suppose first that $A=B$ and $f=\Phi^n_A$; in this case,
we claim that $\fn:=\fm$ will do. Indeed, obviously
$\fn^2=\fn$, and assumption (b) says that $\Phi^n_A(\fm)$
generates the ideal $\fm$; moreover, if $\fn'\subset A$
is another ideal fulfilling these conditions, we see that
$x^p\in\fm$ for every $x\in\fn'$, whence $\fn'\subset\fm$,
since $\fn'$ fulfills condition $(\bB)$, and conversely,
$\fm=\Phi^n_A(\fn')\cdot A\subset\fn'$.

In the general case, we claim that $\fn:=g(\fm)\cdot B$
will do. Indeed, clearly $\fn^2=\fn$ and
$f(\fn)\cdot A=\Phi^n_A(\fm)\cdot A=\fm$, due to assumption
(b); it is also easily seen that $\fn$ fulfills condition
$(\bB)$. Moreover, if $\fn'\subset B$ is another ideal
fulfilling these conditions, it follows that
$$
\Phi^n_B(\fn')\cdot B=g(f(\fn')\cdot A)\cdot B=
g(\fm)\cdot B=\fn
$$
whence $\fn=\fn'$, by the foregoing case.
\end{proof}

\sset\subsubsection{}\label{subsec_old-cartesian-diagr}
Consider now a cartesian diagram of rings
$$
\xymatrix{ A \ar[r]^-{p_1} \ar[d]_{p_2} & A_1 \ar[d]^{\pi_1} \\
A_2 \ar[r]^-{\pi_2} & A_3
}$$
such that $\pi_2$ (and hence $p_1$) is surjective.

\begin{proposition}\label{prop_lift-bas-setup-cartesian}
{\em(i)}\ \ 
In the situation of \eqref{subsec_old-cartesian-diagr}, let
$\fm_1\subset A_1$, $\fm_2\subset A_2$ be two ideals with
$$
\fm_1^2=\fm_1
\qquad
\fm^2_2=\fm_2
\qquad
\fm_1A_3=\fm_2A_3.
$$
\qquad\quad Then there exists a unique ideal $\fm\subset A$
such that
\set\begin{equation}\label{eq_char-smallest}
\fm^2=\fm
\qquad
\fm A_1=\fm_1
\qquad
\fm A_2=\fm_2.
\end{equation}
\begin{enumerate}
\addenu
\item
Moreover, $\fm$ is the smallest of the ideals $J\subset A$
such that $JA_1=\fm_1$ and $JA_2=\fm_2$.
\item
Furthermore, $\fm$ fulfills condition $(\bB)$ if and only
if the same holds for both $\fm_1$ and $\fm_2$.
\end{enumerate}
\end{proposition}
\begin{proof}(i): We shall regard $A$ as a subset of
$A_1\times A_2$, in the natural fashion. Now, set
$I:=\fm_1\times_{A_3}\fm_2$; then $I$ is an ideal
of $A$, and we remark :

\begin{claim}\label{cl_242}
$IA_1=\fm_1$ and $IA_2=\fm_2$.
\end{claim}
\begin{pfclaim} Since $\pi_2$ is surjective, it is easily
seen that $p_1$ restricts to a surjection $I\to\fm_1$.
Next, let $a\in\fm_2$ be any element; we may find an
integer $r\in\N$ and elements $x_1,\dots,x_r\in\fm_1$
and $y_1,\dots,y_r\in A_3$ such that
$a=\sum_{i=1}^ry_i\cdot\pi_1(x_i)$. For every $i=1,\dots,r$,
pick $b_i\in\fm_2$ and $c_i\in A_2$ such that
$\pi_2(b_i)=\pi_1(x_i)$ and $\pi_2(c_i)=y_i$; then
$(x_i,b_i)\in I$ for $i=1,\dots,r$, and
$a-\sum_{i=1}^rp_2(x_i,b_i)\cdot c_i\in\Ker\,\pi_2$.
However, $0\times\Ker\,\pi_2\subset I$, whence $IA_2=\fm_2$.
\end{pfclaim}

\begin{claim}\label{cl_243}
Let $J\subset A$ be any ideal such that $JA_1=\fm_1$
and $JA_2=\fm_2$. Then
$$
(A_1\times_{A_3}\fm_2)\cdot(\fm_1\times_{A_3}A_2)\subset J.
$$
\end{claim}
\begin{pfclaim} Indeed, let $(a_1,a_2),(b_1,b_2)\in A$
be two elements such that $a_2\in\fm_2$ and $b_1\in\fm_1$.
Since $p_1$ is surjective, we may find $(b_1,x)\in J$
such that $\pi_2(x)=\pi_1(b_1)$, and we have
$$
(a_1,a_2)\cdot(b_1,b_2)-(a_1,a_2)\cdot(b_1,x)=
(a_1,a_2)\cdot(0,b_2-x)=(0,a_2\cdot y)
$$
where $y:=b_2-x\in\Ker\,\pi_2$. We are thus reduced to
checking that $(0,ab)\in J$ for every $a\in\fm_2$ and
$b\in\Ker\,\pi_2$. However, say that $a=\sum_{i=1}^r\pi_2(c_i)d_i$
for elements $c_1,\dots,c_r\in J$ and $d_1,\dots,d_r\in A_2$;
it follows that $(0,d_ib)\in A$ for every $i=1,\dots,r$,
and $\sum_{i=1}^rc_i\cdot(0,d_ib)=(0,ab)$, whence the contention.
\end{pfclaim}

From claim \ref{cl_242} we deduce that $I^3A_1=\fm_1^3=\fm_1$,
and likewise, $I^3A_2=\fm_2$; from claim \ref{cl_243}, it
then follows that
$(\fm_1\times_{A_3}A_2)\cdot(A_1\times_{A_3}\fm_2)\subset I^3$.
Especially, $I^2\subset I^3$, so that $I^3=I^2$, and therefore
the ideal $\fm:=I^2$ fulfills conditions \eqref{eq_char-smallest}.

(ii): For every ideal $J$ as in claim \ref{cl_243} we also
get : $\fm\subset J\subset\fm_1\times_{A_3}\fm_2$; if we have
as well $J^2=J$, it follows that $J\subset\fm$, {\em i.e.}
$J=\fm$, as required.

(iii): Due to \cite[Claim 2.19]{Ga-Ra}, it suffices to show,
for every prime $p$, that the following conditions are
equivalent :
\begin{enumerate}
\alphaenu
\item
The $A$-module $\fm/p\fm$ is generated by the $p$-th powers
of its elements
\item
For $i=1,2$, the $A_i$-module $\fm_i/p\fm_i$ is generated
by the $p$-th powers of its elements.
\end{enumerate}
However, it is easily seen that (a)$\Rightarrow$(b). For
the converse, consider the ideal $J\subset A$ generated
by the system $\{x^p~|~x\in\fm\}\cup\{px~|~x\in\fm\}$;
since $p_1$ is surjective, (b) implies that $JA_1=\fm_1$.
Likewise, $JA_2=\fm_2$ : indeed, since $\fm_2/p\fm_2$ is
generated by the $p$-th powers of its elements, it suffices
to check that for every $a\in\fm_2$ we have $a^p\in JA_2$.
But say that $a=\sum_{j=1}^nx_ia_i$ for some
$x_1,\dots,x_n\in\fm$ and $a_1,\dots,a_n\in A_2$; then
$a^p=\sum_{i=1}x_i^pa_i^p+py$ for some $y\in\fm_2$, whence
the contention. In view of (ii), we conclude that $J=\fm$,
whence (a).
\end{proof}

\sset\subsubsection{}
Henceforth, we fix a prime number $p$, and the notation $W(A)$
and $W_{n+1}(A)$ will refer to the ring of {\em $p$-typical}
Witt vector of section \eqref{sec_Witt-Fontaine}, for every
ring $A$ and every integer $n\in\N$. Notice that the functors
$W$ and $W_{n+1}$ introduced in section \ref{sec_Witt-Fontaine}
are defined on the category of topological rings; however, in
this section we shall be interested only in the underlying
rings, and the topologies will play no role; if one wishes,
one may assume that all the rings in this section carry
the discrete topology. We consider first the case where
$p$ is nilpotent on $A$ :

\begin{corollary}\label{cor_alm-str-on-Witt}
Let $A$ be any ring such that $p^kA=0$ for
some $k\in\N$, and $\fm\subset A$ an ideal with
$\fm^2=\fm$. For every $n\in\N$ we have :
\begin{enumerate}
\item
There exists a unique ideal $\fn_{n+1}\subset W_{n+1}A$ such that
$$
\fn_{n+1}^2=\fn_{n+1}
\qquad\text{and}\qquad
\bar\bomega_0(\fn_{n+1})=\fm.
$$
\item
The image of $\fn_{n+2}$ under the projection $W_{n+2}A\to W_{n+1}A$
agrees with $\fn_{n+1}$.
\end{enumerate}
\end{corollary}
\begin{proof} (See remark \ref{rem_Witt-limit}(i) for the
definition of the ring homomorphism $\bar\bomega_0$.)

(i) follows immediately from corollary \ref{cor_int-with-nilker}
and lemma \ref{lem_lift-a-struct-over-nilp}.

(ii) follows immediately from (i).
\end{proof}

For a general ring $A$, we may state :

\begin{proposition}\label{prop_alm-str-on-Witt}
Let $A$ be a ring, $n\in\N$ an integer, $\fm\subset A$ an
ideal with $\fm=\fm^2$ and fulfilling condition $(\bB)$
of\/ \cite[\S2.1.6]{Ga-Ra}. Set
$\fn_{n+1}=\{\underline a\in W_{n+1}A~|~a_0,\dots,a_n\in\fm\}$.
We have :
\begin{enumerate}
\item
$\fn_{n+1}$ is an ideal of\/ $W_{n+1}A$ with $\fn_{n+1}^2=\fn_{n+1}$
and fulfilling condition $(\bB)$.
\item
Moreover, $\fn_{n+1}$ is the unique ideal with $\fn_{n+1}^2=\fn_{n+1}$
and such that
\set\begin{equation}\label{eq_ghost-used-twice}
\bar\bomega_i(\fn_{n+1})\cdot A=\fm
\qquad
\text{for $i=0,\dots,n$}.
\end{equation}
\end{enumerate}
\end{proposition}
\begin{proof} We consider the ring homomorphism
$$
\pi_n:W_{n+1}A\to A^{n+1}
$$
of proposition \ref{prop_int-with-nilker}(i). By proposition
\ref{prop_int-with-nilker}(i), the kernel of the surjection
$W_{n+1}A\to R:=\Img\,\pi_n$ is nilpotent, hence every basic
setup $(R,\fm_R)$ lifts uniquely to a basic setup
$(W_{n+1}A,\tilde\fm)$, and $\tilde\fm$ satisfies condition
$(\bB)$ if and only if the same holds for $\fm_R$ (lemma
\ref{lem_lift-a-struct-over-nilp}(i,iii)). Moreover, by
proposition \ref{prop_int-with-nilker}(ii), the subring
$R$ of $A^{n+1}$ contains the ideal $p^{n+1}A^{n+1}$, and lies
in the larger subring $S\subset A^{n+1}$ consisting of all
sequences $(a_0,\dots,a_n)$ such that
$a_i\equiv a_0^{p^i}\pmod{pA}$ for $i=0,\dots,n$; there
follows a cartesian diagram of rings
\set\begin{equation}\label{eq_another-cart-structures}
{\diagram R \ar[r] \ar[d] & R_0:=R/p^{n+1}A^{n+1} \ar[d] \\
S \ar[r] & S_0:=S/p^{n+1}A^{n+1}
\enddiagram}
\end{equation}
whose horizontal arrows are surjections. By proposition
\ref{prop_lift-bas-setup-cartesian}(i,iii), the datum
of $(R,\fm_R)$ is then equivalent to that of a pair of
basic setups $(R_0,\fm_{R_0})$ and $(S,\fm_S)$ with
$\fm_{R_0}S_0=\fm_SS_0$. Also, $\fm_R$ fulfills condition
$(\bB)$ if and only if the same holds for both $\fm_{R_0}$
and $\fm_S$. Next, set $A_0:=A/pA$, and notice that
we have a cartesian diagram of rings
\set\begin{equation}\label{eq_cartesian-structures}
{\diagram S \ar[r] \ar[d]_\psi & A^{n+1} \ar[d] \\
A_0 \ar[r]^-\phi & A_0^{n+1} 
\enddiagram}
\end{equation}
where $\phi$ is given by the rule :
$x\mapsto(x,x^p,\dots,x^{p^n})$ for every $x\in A_0$;
the top horizontal arrow is the inclusion map, and the
right vertical arrow is the projection. Moreover, $\psi$
is a surjection. Hence -- again by virtue of proposition
\ref{prop_lift-bas-setup-cartesian}(i) -- the datum of
$(S,\fm_S)$ is equivalent to that of a system of basic
setups $(A,\fm_i)$ such that
$$
\fm_iA_0=\Phi^i_{A_0}(\fm_0A_0)\cdot A_0
\qquad
\text{for $i=1,\dots,n$}
$$
where $\Phi_{A_0}:A_0\to A_0$ is the Frobenius endomorphism.
Especially, if $(A,\fm)$ is a basic setup on $A$ such that
$\fm$ fulfills condition $(\bB)$, then we can take
$\fm_i:=\fm$ for every $i=0,\dots,n$, and in this way we
obtain a natural basic setup $(S,\fm_S)$ on $S$ with
$\fm_SA^{n+1}=\fm^{n+1}$, and such that condition $(\bB)$
still holds for $\fm_S$. Furthermore, the restriction
of $\psi$ to $R$ factors through a surjection $R_0\to A_0$
whose kernel is also nilpotent, so that the basic setup
$(A_0,\fm A_0)$ lifts uniquely to a basic setup $(R_0,\fm_{R_0})$
with $\fm_{R_0}$ fulfilling condition $(\bB)$ (lemma
\ref{lem_lift-a-struct-over-nilp}(i)). Lastly, both ideals
$\fm_{R_0}S_0$ and $\fm_SS_0$ lift $\fm A_0$ along the map
$S_0\to A_0$ induced by $\psi$, whose kernel is also nipotent;
by invoking again \ref{lem_lift-a-struct-over-nilp}(i), we
deduce that $\fm_{R_0}S_0=\fm_SS_0$. Summing up, we have finally
associated with the basic setup $(A,\fm)$ a natural basic
setup $(W_{n+1}A,\tilde\fm)$ which is the unique one with
property \eqref{eq_ghost-used-twice}, and with $\tilde\fm$
fulfilling condition $(\bB)$. To conclude the proof of both
(i) and (ii), it then suffices to check that $\tilde\fm=\fn_{n+1}$.
However, it is clear that $\fn_{n+1}$ is an ideal of $W_{n+1}A$
enjoying property \eqref{eq_ghost-used-twice}, since
$\bar\bomega_i(\fn_{n+1})$ is contained in $\fm$ and contains
the system $(x^{p^i}~|~x\in\fm)$, which generates $\fm$, under
condition $(\bB)$. It remains to check that $\fn_{n+1}=\fn_{n+1}^2$.
Now, let us define a descending filtration by subideals of
$\fn_{n+1}$, by setting $\Fil_i\fn_{n+1}:=\fn_{n+1}\cap\bar V_iA$,
where $\bar V_iA:=V_iA/V_{n+1}A$ for every $i=0,\dots,n$. Let
also $\gr_\bullet\fn_{n+1}$ be the associated graded $W_{n+1}$-module.
We are then further reduced to showing that
$\fn_{n+1}\cdot\gr_i\fn_{n+1}=\gr_i\fn_{n+1}$ for $i=0,\dots,n$.
However, we have a natural identification of $W_{n+1}A$-modules :
$\gr_i\fn_{n+1}\isom\fm$, for the $W_{n+1}A$-module structure on
$\fm$ induced by restriction of scalars along the ghost map
$\bar\bomega_i:W_{n+1}A\to A$ (claim \ref{cl_W_n-structure}).
Since $\fm^2=\fm$, we then come down to the assertion that
$\bar\bomega_i(\fn_{n+1})\cdot A=\fm$, which was already
remarked.
\end{proof}

\sset\subsubsection{}\label{subsec_almost-Witts}
Let $(A,\fm)$ be a basic setup such that either $p^kA=0$
for some integer $k\in\N$, or else such that $\fm$ fulfills
condition $(\bB)$. For every $n\in \N$ let
$\fn_{n+1}\subset W_{n+1}A$ be the ideal provided by
corollary \ref{cor_alm-str-on-Witt}(i), or respectively
by proposition \ref{prop_alm-str-on-Witt}; especially,
$(W_{n+1}A,\fn_{n+1})$ is also a basic setup. Let also
$f:B\to C$ be a morphism of $A$-algebras such that
$f^a:B^a\to C^a$ is an isomorphism of $(A,\fm)^a$-algebras.
Then we claim that $f$ induces an isomorphism of
$(W_{n+1}A,\fn_{n+1})^a$-algebras
$$
(W_{n+1}f)^a:(W_{n+1}B)^a\isom(W_{n+1}C)^a
\qquad
\text{for every $n\in\N$}.
$$
Indeed, arguing by induction on $n\in\N$, we are easily
reduced to checking that $f$ induces an isomorphism of
$(W_{n+1}A)^a$-modules $(V_nB/V_{n+1}B)^a\isom(V_nC/V_{n+1}C)^a$ for
every $n\in\N$. However, in light of claim \ref{cl_W_n-structure}
we have a commutative diagram of $W_{n+1}A$-modules
$$
\xymatrix{ B \ar[r] \ar[d]_f & V_nB/V_{n+1}B \ar[d] \\
C \ar[r] & V_nC/V_{n+1}C
}$$
whose horizontal arrows are isomorphisms of $W_{n+1}A$-modules,
and where the $W_{n+1}A$-module structures on $B$ and $C$ are
induced by the $n$-th ghost maps. By assumption,
$\fm$ annihilates $\Ker\,f$ and $\Coker\,f$; since
$\bar\bomega_n(\fn_{n+1})\subset\fm$, it follows that
$\fn_{n+1}$ annihilates these $W_{n+1}A$-modules as well,
whence the contention. Thus, we obtain a well defined
functor
$$
W_{n+1}:(A,\fm)^a\Alg\to(W_{n+1}A,\fn_{n+1})^a\Alg
\qquad
B^a\mapsto(W_{n+1}B)^a
\qquad
\text{for every $n\in\N$}.
$$
Moreover, for every $k\leq n$, the ghost component
$\bar\bomega_k:W_nB\to B$ induces pull-back functors
$$
\bar\bomega^*_k:B^a\Alg\to W_{n+1}B^a\Alg
\qquad
\bar\bomega^*_k:B^a\Mod\to W_{n+1}B^a\Mod.
$$
Namely, for every $B$-algebra $C$ we let $\bar\bomega_k^*C$
be the $W_{n+1}B$-algebra whose underlying ring is $C$,
and whose structure morphism is the composition of
$\bar\bomega_k$ and the structure morphism $B\to C$
of $C$. Since $\bar\bomega_k(\fn_{n+1})\subset\fm$, it
is clear that this functor on $B$-algebras descends to
a well defined functor $\bar\bomega^*_k$ on $B^a$-algebras,
as stated. Likewise one argues to define the functor
$\bar\bomega^*_k$ on $B^a$-modules. Furthermore, we have
a well defined ideal
$$
\bar V_kB^a:=(V_kB/V_{n+1}B)^a\subset W_{n+1}B^a
\qquad
\text{for every $k=0,\dots,n$}
$$
and claim \ref{cl_W_n-structure} yields a natural
isomorphism of $W_{n+1}B^a$-modules :
\set\begin{equation}\label{eq_from-cl-W_n-structure}
\bar V_kB^a/\bar V_{k+1}B^a\isom\bar\bomega_k^*B^a
\qquad
\text{for every $k=0,\dots,n$}.
\end{equation}

\sset\subsubsection{}\label{subsec_compute-al.els-of-Witts}
In the situation of \eqref{subsec_almost-Witts}, suppose that
$\fm$ fulfills condition $(\bB)$, and let $B$ be any $A$-algebra;
denote by
$$
B^a_{!!}\xrightarrow{\eps_B}B\xrightarrow{\eta_B}B^a_*
\qquad\text{and}\qquad
(W_{n+1}B)^a_{!!}\xrightarrow{\eps_{W_{n+1}B}}W_{n+1}B
\xrightarrow{\eta_{W_{n+1}B}}(W_{n+1}B)^a_*
$$
the units and counits of adjunction. We have :

\begin{proposition}
In the situation of \eqref{subsec_compute-al.els-of-Witts},
the following holds :
\begin{enumerate}
\item
There exist isomorphisms of\/ $W_{n+1}A$-algebras
$$
\omega:W_{n+1}(B^a_*)\isom(W_{n+1}B)^a_*
\qquad\text{and}\qquad
\tau:W_{n+1}(B^a_{!!})\isom(W_{n+1}B)^a_{!!}
$$
such that $\omega\circ W_{n+1}(\eta_B)=\eta_{W_{n+1}(B)}$ and
$\eps_{W_{n+1}(B)}\circ\tau=W_{n+1}(\eps_B)$.
\item
Moreover, we have a natural isomorphism of non-unital rings
(see remark {\em\ref{rem_non-unital-ring}})
$$
\tilde\fn_{n+1}:=\fn_{n+1}\otimes_{W_{n+1}A}\fn_{n+1}\isom
W_{n+1}(\tilde\fm).
$$
\end{enumerate}
\end{proposition}
\begin{proof}(i): We consider the commutative diagrams
$$
\xymatrix@C+35pt{
W_{n+1}B \ar[r]^-{W_{n+1}(\eta_B)} \ar[d]_{\eta_{W_{n+1}B}} &
W_{n+1}(B^a_*) \ar[d]^{\eta_{W_{n+1}(B^a_*)}} &
W_{n+1}(B^a_{!!})^a_{!!} \ar[r]^-{(W_{n+1}\eps_B)^a_{!!}}
\ar[d]_{\eps_{W_{n+1}(B^a_{!!})}} &
W_{n+1}(B)^a_{!!} \ar[d]^{\eps_{W_{n+1}(B)}} \\
(W_{n+1}B)^a_* \ar[r]^-{(W_{n+1}\eta_B)^a_*} & W_{n+1}(B^a_*)^a_* &
W_{n+1}(B^a_{!!}) \ar[r]^-{W_{n+1}(\eps_B)} & W_{n+1}(B)
}$$
and notice that $(W_{n+1}\eta_B)^a_*$ and $(W_{n+1}\eps_B)^a_{!!}$
are isomorphisms, by the discussion of \eqref{subsec_almost-Witts}.
Hence, it suffices to show that $\eta_{W_{n+1}(B^a_*)}$ and
$\eps_{W_{n+1}(B^a_{!!})}$ are isomorphisms, which follows from :

\begin{claim} If $\eta_B$ (resp. $\eps_B$) is an isomorphism,
the same holds for $\eta_{W_{n+1}B}$ (resp. $\eps_{W_{n+1}B}$).
\end{claim}
\begin{pfclaim} Let $\gr_\bullet(W_{n+1}B)$ be the graded ring
associated with the filtration $(\bar V_kB~|~k=0,\dots,n)$
of $W_{n+1}B$; consider also the filtration
$((\bar V_kB^a)_*~|~k=0,\dots,n)$ of $(W_{n+1}B)^a_*$, and
let $\gr_\bullet(W_{n+1}B)^a_*$ be the associated graded ring.
We get a commutative diagram
$$
\xymatrix{ & \gr_i(W_{n+1}B)
\ar[dl]_{\gr_i(\eta_{W_{n+1}B})} \ar[dr]^{\eta_{\gr_i(W_{n+1}B)}} \\
\gr_i(W_{n+1}B)^a_* \ar[rr] & & (\gr_iW_{n+1}B)^a_*
}$$
whose bottom horizontal arrow is injective, since the functor
$(-)^a_*$ is left exact. Here $\eta_{\gr_i(W_{n+1}B)}$ is again
the unit of adjunction, which, by virtue of
\eqref{eq_from-cl-W_n-structure}, is naturally identified with
the unit of adjunction
$$
\eta_{\bar\bomega_i^*B}:\bar\bomega_i^*B\to(\bar\bomega_i^*B)^a_*
\qquad
\text{for every $i=0,\dots,n$}.
$$
Let now $\cB$ be the category of basic setups, and $\cB\Alg\to\cB$
(resp. $\cB^a\Alg\to\cB$) the fibred and cofibred category of
$\cB$-algebras (resp. of almost $\cB$-algebras), as in
\cite[\S3.5]{Ga-Ra}. The localization $\cB\Alg\to\cB^a\Alg$
admits a right adjoint
$$
(-)_*:\cB^a\Alg\to\cB\Alg
\qquad
((V,\fm),R)\mapsto((V,\fm),R_*)
$$
which is a $\cB$-cartesian functor (\cite[\S3.5.4]{Ga-Ra}).
By \eqref{eq_ghost-used-twice}, we have the morphism
$\bar\bomega_i:(W_{n+1}A,\fn_{n+1})\to(A,\fm)$ in $\cB$, whence
a commutative diagram in $\cB\Alg$ :
$$
\xymatrix{ ((W_{n+1}A,\fn_{n+1}),\bar\bomega^*_iB) \ar[r]
\ar[d]_{\eta_{\bar\bomega^*_iB}} & ((A,\fm),B) \ar[d]^{\eta_B} \\
((W_{n+1}A,\fn_{n+1}),(\bar\bomega^*_iB)^a_*) \ar[r] & ((A,\fm),B^a_*)
}$$
whose top horizontal arrow is given by the identity map of
$\bar\bomega^*_iB$, and is therefore a cartesian morphism of
$\cB\Alg$. Then also the bottom vertical arrow is a cartesian
morphism, since $(-)_*$ is a cartesian functor. We conclude
that $\eta_{\gr_i(W_{n+1}B)}$ is an isomorphism, under our assumptions.
Then the same holds for $\gr_i(\eta_{W_{n+1}B})$, for every
$i=0,\dots,n$, whence the assertion for $\eta_{W_{n+1}B}$.

To show the assertion concerning $\eps_{W_{n+1}B}$, we argue by
induction on $n\in\N$; the case $n=0$ is trivial. Let then $n>0$,
and suppose that the assertion is already known for $\eps_{W_nB}$;
recall that for every $i\in\N$, the ring $(W_{i+1}B)^a_{!!}$ is the
cokernel of a map of $W_{i+1}A$-modules
$$
\delta_i:\tilde\fn_{i+1}\to
W_{i+1}A\oplus(\tilde\fn_{i+1}\otimes_{W_{i+1}A}W_{i+1}B)
\qquad
x_\bullet\otimes y_\bullet\mapsto
(x_\bullet\cdot y_\bullet,x_\bullet\otimes y_\bullet\otimes 1).
$$
We consider the commutative ladder with right exact rows :
$$
\xymatrix{
\tilde\fn_{n+1} \ar[r]^-{\delta_n} \ar[d]_\alpha & W_{n+1}A\oplus
(\tilde\fn_{n+1}\otimes_{W_{n+1}A}W_{n+1}B)
\ar[r] \ar[d]^\beta & (W_{n+1}B)^a_{!!} \ar[r] \ar[d]^\gamma & 0 \\
\tilde\fn_n \ar[r]^-{\delta_{n-1}} &
W_nA\oplus(\tilde\fn_n\otimes_{W_nA}W_nB) \ar[r] &
(W_nB)^a_{!!} \ar[r] & 0
}$$
whose vertical arrows are induced by the projections
$\pi^A_n:W_{n+1}A\to W_nA$ and $\pi^B_n:W_{n+1}B\to W_nB$ (notice
that the almost structure used to compute $(W_{n+1}B)^a_{!!}$ is
the one relative to the basic setup $(W_{n+1}A,\fn_{n+1})$, whereas
for $(W_nB)^a_{!!}$ we use the almost structure relative to the
basic setup $(W_nA,\fn_n)$). Since $\pi^A_n(\fn_{n+1})=\fn_n$,
by \cite[Rem.2.1.4(ii)]{Ga-Ra} we have a natural isomorphism
$$
\tilde\fn_{n+1}\otimes_{W_{n+1}A}W_nA\isom\tilde\fn_n
\qquad
x_\bullet\otimes y_\bullet\otimes a_\bullet\mapsto
a_\bullet\cdot\pi^A_n(x_\bullet)\otimes\pi^A_n(y_\bullet)
$$
which identifies $\alpha$ with
$\tilde\fn_{n+1}\otimes_{W_{n+1}A}\pi^A_n$, and $\beta$ with
$\pi^A_n\oplus(\tilde\fn_{n+1}\otimes_{W_{n+1}A}\pi^B_n)$.
Taking claim \ref{cl_W_n-structure} into account, we deduce
natural $W_{n+1}A$-linear surjections :
\set\begin{equation}\label{eq_one-more}
\tilde\fn_{n+1}\otimes_{W_{n+1}A}\bar\bomega_n^*(A)\to\Ker\,\alpha
\qquad
\bar\bomega_n^*(A)\oplus
(\tilde\fn_{n+1}\otimes_{W_{n+1}A}\bar\bomega^*_n(B))\to\Ker\,\beta.
\end{equation}
Then, in light of \eqref{eq_ghost-used-twice} and
\cite[Rem.2.1.4(ii)]{Ga-Ra}, we have as well the natural isomorphism
$$
\tilde\fn_{n+1}\otimes_{W_{n+1}A}\bar\bomega_n^*(A)\isom\tilde\fm
\qquad
x_\bullet\otimes y_\bullet\otimes a\mapsto
a\cdot\bar\bomega_n(x_\bullet)\otimes\bar\bomega_n(y_\bullet)
$$
which identifies \eqref{eq_one-more} with $W_{n+1}A$-linear
surjections :
\set\begin{equation}\label{eq_ker-beta}
\bar\bomega_n^*(\tilde\fm)\to\Ker\,\alpha
\qquad
\bar\bomega^*_n(A\oplus(\tilde\fm\otimes_AB))\to\Ker\,\beta.
\end{equation}
By a direct inspection, we then get a commutative diagram :
$$
\xymatrix@C+20pt{ \bar\bomega_n^*(\tilde\fm)
\ar[r]^-{\bar\bomega_n^*(\delta_0)} \ar[d] &
\bar\bomega^*_n(A\oplus(\tilde\fm\otimes_AB)) \ar[d] \\
\Ker\,\alpha \ar[r]^-{\delta'_n} & \Ker\,\beta
}$$
whose vertical arrows are the surjections \eqref{eq_ker-beta},
and $\delta'_n$ is the restriction of $\delta_n$. Since $\alpha$
and $\beta$ are surjective, the induced map
$\Coker\,\delta'_n\to\Ker\,\gamma$ is bijective, by the snake
lemma, and $\gamma$ is surjective, so finally we get a right exact
sequence of $W_{n+1}A$-modules :
$$
\bar\bomega_n^*(B^a_{!!})\xrightarrow{\rho}(W_{n+1}B)^a_{!!}
\xrightarrow{\gamma}(W_nB)^a_{!!}\to 0.
$$
We consider then the diagram of $W_{n+1}A$-modules :
$$
\xymatrix{ &
\bar\bomega_n^*(B^a_{!!}) \ar[r]^-\rho \ar[d]_{\bar\bomega_n^*\eps_B} &
(W_{n+1}B)^a_{!!} \ar[r]^-\gamma \ar[d]^{\eps_{W_{n+1}B}} &
(W_nB)^a_{!!} \ar[r] \ar[d]^{\eps_{W_nB}} & 0 \\
0 \ar[r] & \bar\bomega_n^*B \ar[r]^-{\rho'} &
W_{n+1}B \ar[r]^-{\gamma'} & W_nB \ar[r] & 0
}$$
where $\gamma'$ is the projection, and $\rho'$ is the
natural identification of $\bar\bomega_n^*B$ with $\Ker\,\gamma'$
given by claim \ref{cl_W_n-structure}. Thus, the bottom horizontal
sequence is exact, and by inductive assumption the first and
third vertical arrows are bijective; therefore, in order to
conclude the proof, it will suffice to check the diagram
commutes. The commutativity of the right square subdiagram
is clear. To see the commutativity of the left square subdiagram,
let $w\in\bar\bomega_n^*(B^a_{!!})$ be any element; hence $w$ is
the class $[a,z]$ of a pair $(a,z)\in A\oplus(\tilde\fm\otimes_AB)$,
and we may assume that $z=x\otimes y\otimes b$ for some $x,y\in\fm$
and $b\in B$. Then, by inspecting the constructions, we see that
$w':=\rho(w)$ is the class $[(0,\dots,0,a),(0,\dots,0,xyb)]$ of
the corresponding pair in
$W_{n+1}A\oplus(\tilde\fn_{n+1}\otimes_{W_{n+1}A}W_{n+1}B)$. On the
other hand, $\bar\bomega_n^*\eps_B(w)=a+xyb$; then, again by
claim \ref{cl_W_n-structure}, we deduce that
$\eps_{W_{n+1}B}(w')=(0,\dots,0,a+xyb)=\rho'(a+xyb)$, as required.
\end{pfclaim}

(ii): For every $n\in\N$, we have a natural map of
non-unital $W_{n+1}A$-algebras :
$$
\mu_n:W_{n+1}(\fm)\otimes_{W_{n+1}A}W_{n+1}(\fm)\to W_{n+1}(\tilde\fm)
\qquad
a_\bullet\otimes b_\bullet\mapsto
(P_i(a_\bullet\otimes 1,1\otimes b_\bullet)~|~i=0,\dots,n)
$$
for the non-unital ring structure of $W_{n+1}(\fm)=\fn_{n+1}$ and
$W_{n+1}(\tilde\fm)$ induced by those of $\fm$ and respectively
$\tilde\fm$ (remark \ref{rem_non-unital-ring}). We show by induction
on $n\in\N$, that $\mu_n$ is an isomorphism; the assertion is trivial
for $n=0$, since $\mu_0=\one_{\tilde\fm}$. Thus, let $n>0$, and suppose
that the assertion is already known for $n-1$; we consider the
commutative diagram of $W_{n+1}A$-modules :
$$
\xymatrix{
\fn_{n+1}\otimes_{W_{n+1}A}\fn_{n+1} \ar[r]^-{\mu_n} \ar[d]_\alpha &
W_{n+1}(\tilde\fm) \ar[d]^\beta \\
\fn_n\otimes_{W_nA}\fn_n \ar[r]^-{\mu_{n-1}} & W_n(\tilde\fm)
}$$
whose left vertical arrows is induced by the projection
$W_{n+1}(\fm)\to W_n(\fm)$, and whose right vertical arrow
is likewise the natural projection. Claim \ref{cl_W_n-structure}
yields a natural isomorphism of $W_{n+1}A$-modules :
\set\begin{equation}\label{eq_etc}
\bar\bomega_n^*(\tilde\fm)\isom\Ker\,\beta.
\end{equation}

\begin{claim}\label{cl_done}
We have a commutative diagram of $W_{n+1}A$-modules :
$$
\xymatrix@C+20pt{ \bar\bomega_n^*(\tilde\fm)
\ar[r]^-{\bar\bomega_n^*(\mu_0)} \ar[d] &
\bar\bomega_n^*(\tilde\fm) \ar[d] \\
\Ker\,\alpha \ar[r] & \Ker\,\beta
}$$
whose bottom horizontal arrow is the restriction of $\mu_n$,
and whose right (resp. left) vertical arrow is \eqref{eq_etc}
(resp. is \eqref{eq_ker-beta}).
\end{claim}
\begin{pfclaim} By construction, the left vertical arrow is
characterized as the map such that
$a\cdot\bar\bomega_n(x_\bullet)\otimes\bar\bomega_n(y_\bullet)
\mapsto(0,\dots,0,a)\cdot(x_\bullet\otimes y_\bullet)$ for
every $a\in A$ and every $x_\bullet,y_\bullet\in\fn_{n+1}$.
On the other hand, the right vertical arrow is the map
such that $x\otimes y\mapsto(0,\dots,0,x\otimes y)$ for
every $x,y\in\fm$. Hence, we need to check the identity :
$$
\mu_n((0,\dots,0,a)\cdot(x_\bullet\otimes y_\bullet))=
(0,\dots,0,a\cdot\bar\bomega_n(x_\bullet)\otimes\bar\bomega_n(y_\bullet))
$$
for every $a\in A$ and every $x_\bullet,y_\bullet\in\fn_{n+1}$.
To this aim, endow $R:=A\oplus\fm$ with its natural $A$-algebra
structure as in remark \ref{rem_non-unital-ring}(iii), and
recall that $\mu_n$ is the restriction of the map of unital
$W_{n+1}A$-algebras
$$
W_{n+1}R\otimes_{W_{n+1}A}W_{n+1}R\to W_{n+1}(R\otimes_AR)
\qquad
x_\bullet\otimes y_\bullet\mapsto
(x_\bullet\otimes 1)\cdot(1\otimes y_\bullet).
$$
Thus, we are reduced to checking that
$$
((V^n_R(\tau_R(a))\cdot x_\bullet)\otimes 1)\cdot(1\otimes y_\bullet)=
V^n_{R\otimes_AR}(\tau_{R\otimes_AR}(a\cdot\bar\bomega_n(x_\bullet)\otimes
\bar\bomega_n(y_\bullet)))
$$
for every $a\in A$ and every $x_\bullet,y_\bullet\in W_{n+1}R$, where
$V_R$ and $\tau_R$ denote respectively the Verschiebung and the
Teichm\"uller maps, and likewise for $V_{R\otimes_AR}$ and
$\tau_{R\otimes_AR}$ (see \eqref{subsec_V-and-F}). Let also $F_R$
be the Frobenius map; using proposition \ref{prop_V_A-and-F_A}(ii)
we compute :
$$
\begin{aligned}
((V^n_R(\tau_R(a))\cdot x_\bullet)\otimes 1)\cdot(1\otimes y_\bullet)&=
(V^n_R(\tau_R(a)\cdot F^n_R(x_\bullet))\otimes 1)\cdot(1\otimes y_\bullet) \\
&= V^n_R((\tau_R(a)\cdot F^n_R(x_\bullet))\otimes 1)\cdot(1\otimes y_\bullet) \\
&= V^n_{R\otimes_AR}(((\tau_R(a)\cdot F^n_R(x_\bullet))\otimes 1)\cdot
(1\otimes F^n_R(y_\bullet))) \\
&= V^n_{R\otimes_AR}(\tau_{R\otimes_AR}(a)\cdot F^n_{R\otimes_AR}(x_\bullet\otimes 1)
\cdot F^n_{R\otimes_AR}(1\otimes y_\bullet)).
\end{aligned}
$$
Hence, we are reduced to checking that :
$$
\bar\bomega_n(x_\bullet)\otimes\bar\bomega_n(y_\bullet)=
\bar\bomega_0(F^n_{R\otimes_AR}(x_\bullet\otimes 1)
\cdot F^n_{R\otimes_AR}(1\otimes y_\bullet))
\qquad
\text{in $W_{n+1}(R\otimes_AR)$}.
$$
The latter follows easily from \eqref{eq_basic-F} : details
left to the reader.
\end{pfclaim}

Claim \ref{cl_done} implies that also \eqref{eq_ker-beta}
is an isomorphism, and the same for $\mu_{n-1}$, by inductive
assumption; so the same follows for $\mu_n$ and the proof is
concluded.
\end{proof}

\sset\subsubsection{}\label{subsec_explain-Frob}
Let now $A_0$ be any $\F_p$-algebra, $(A_0,\fm_0)$ any
basic setup, and $R$ any $A_0$-algebra. Recall that the
Frobenius endomorphism $\Phi_R:R\to R$ induces an endofunctor
$$
\Phi^{k*}_R:R\Alg\to R\Alg
\qquad
\text{for every $k\in\N$}
$$
that assigns to every $R$-algebra $S$ the $R$-algebra
$\Phi^{k*}_RS$ whose underlying ring is the same as $S$, and
whose structure morphism is the composition of $\Phi^k_R$
with the structure morphism $R\to S$ of $S$. Moreover,
the Frobenius endomorphism $\Phi^k_S$ induces a natural
transformation
$$
\Phi^k_{S/R}:S\otimes_R\Phi^{k*}_RR\to\Phi^{k*}_RS
\qquad
x\otimes y\mapsto x^{p^k}y
\qquad
\text{for every $R$-algebra $S$}
$$
Let $R^a$ be the $(A_0,\fm)^a$-algebra represented
by $R$; then $\Phi^{k*}_R$ descends to an endofunctor
$$
\Phi_{R^a}^{k*}:R^a\Alg\to R^a\Alg
$$
that assigns to every $R^a$-algebra $S^a$ the $R^a$-algebra
$\Phi_{R^a}^{k*}(S^a):=(\Phi_R^{k*}S)^a$, where $S$ is any
$R$-algebra representing $S^a$ : see \cite[\S3.5.7]{Ga-Ra}.
Likewise, the Frobenius endomorphism $\Phi^k_S$ induces a
natural transformation
$$
\Phi^k_{S^a/R^a}:=(\Phi^k_{S/R})^a:
(S^a\otimes_{R^a}\Phi^{k*}_{R^a}R^a)\to\Phi^{k*}_{R^a}S^a
\qquad
\text{for every $R^a$-algebra $S^a$}.
$$
Moreover, notice that $\Phi^{k+1}_{S^a/R^a}$ also equals the
composition of
$\Phi_{S^a/R^a}\otimes_{\Phi^*_{R^a}R^a}\Phi^{k+1*}_{R^a}R^a$ :
$$
S^a\otimes_{R^a}\Phi^{k+1*}_{R^a}R^a\isom
S^a\otimes_{R^a}\Phi^*_{R^a}R^a\otimes_{\Phi^*_{R^a}R^a}\Phi^{k+1*}_{R^a}R^a
\to\Phi^*_{R^a}S^a\otimes_{\Phi^*_{R^a}R^a}\Phi^{k+1*}_{R^a}R^a
$$
together with $\Phi^*_{R^a}(\Phi^k_{S^a/R^a}):
\Phi^*_{R^a}S^a\otimes_{\Phi^*_{R^a}R^a}\Phi^{k+1*}_{R^a}R^a=
\Phi^*_{R^a}(S^a\otimes_{R^a}\Phi^{k*}_{R^a}R^a)\to
\Phi^{k+1*}_{R^a}S^a$. Especially, if $\Phi_{S^a/R^a}$ is
an isomorphism, the same holds for $\Phi^k_{S^a/R^a}$,
for every $k\in\N$.

\begin{remark}\label{rem_explain-Frob}
In the situation of \eqref{subsec_explain-Frob},
suppose that $\fm$ fulfills condition $(\bB)$ of
\cite[\S2.1.6]{Ga-Ra}, and let $f:R\to S$ be any
weakly \'etale morphism of $A_0^a$-algebras. Then
$\Phi_{S/R}$ is an isomorphism. Indeed, the assertion
is \cite[Th.3.5.13(ii)]{Ga-Ra} in case $\fm\otimes_{A_0}\fm$
is a flat $A_0$-module, but by direct inspection one sees
easily that the proof of {\em loc.cit.} is valid more
generally, whenever condition $(\bB)$ holds.
\end{remark}

\begin{theorem}\label{th_Witt-and-etale-maps}
In the situation of \eqref{subsec_almost-Witts}, let
$f:B\to C$ be a flat morphism of $A^a$-algebras; set
$B_0:=B/pB$, $C_0:=C/pC$, and suppose that $\Phi_{C_0/B_0}$
is an isomorphism. Then:
\begin{enumerate}
\item
$f$ induces a flat morphism of\/ $W_{n+1}A^a$-algebras
$$
W_{n+1}f:W_{n+1}B\to W_{n+1}C
\qquad
\text{for every $n\in\N$}.
$$
\item
For every $n\in\N$ and every $i=0,\dots,n$, the following
diagram is cocartesian :
$$
\xymatrix{ W_{n+1}B \ar[rr]^-{W_{n+1}f} \ar[d]_{\bar\bomega_i} & &
W_{n+1}C \ar[d]^{\bar\bomega_i} \\
\bar\bomega_i^*B \ar[rr]^-{\bar\bomega_i^*f} & & \bar\bomega_i^*C.
}$$
\end{enumerate}
\end{theorem}
\begin{proof} We prove first both assertions in the case where
$p^kA=0$ for some integer $k\in\N$. For every $j\in\N$, and
every $A^a$-algebra $R$ set
$$
W_{n,j}R:=W_{n+1}(R)/p^j\bar V_n(R)
\qquad
\text{with $\bar V_n(R):=V_nR/V_{n+1}R\subset W_{n+1}R$}.
$$
Hence, $W_{n,0}R=W_nR$, and in light of corollary
\ref{cor_int-with-nilker}, we see that $W_{n,j}R=W_{n+1}R$
for every sufficiently large $j\in\N$. We shall show, by
induction on $j$ and $n$, the following assertions :

$\mathrm{(a)}_{n,j}$\ \
The induced map $W_{n,j}f:W_{n,j}B\to W_{n,j}C$ is flat for
every $j,n\in\N$.

$\mathrm{(b)}_{n,j}$\ \
The morphism $f$ induces a cocartesian diagram of
$W_{n+1}B$-algebras :
$$
\xymatrix{
W_{n,j}B \ar[r] \ar[d]_{W_{n,j}f} &
\bar\bomega^*_0B \ar[d]^{\bar\bomega^*_0f} \\
W_{n,j}C \ar[r] & \bar\bomega^*_0C.
}$$
Both assertions are trivial for $n=0$ and every $j\in\N$.
Let then $n,j\in\N$ be any integers, and suppose that
$\mathrm{(a)}_{n,j}$ and $\mathrm{(b)}_{n,j}$ hold. According
to \eqref{eq_from-cl-W_n-structure}, the kernel $K_{n,j}$
of the projection $W_{n,j+1}B\to W_{n,j}B$ is isomorphic,
as a $W_{n+1}B$-module, to $\bar\bomega^*_n(p^jB/p^{j+1}B)$.
Likewise for the kernel $K'_{n,j}$ of the projection
$W_{n,j+1}C\to W_{n,j}C$, and the induced morphism
$K_{n,j}\to K'_{n,j}$ corresponds to the restriction of
$\bar\bomega^*_n(f\otimes_\Z\Z/p^{j+1}\Z)$, under these
identifications. By the (almost version of the) local
flatness criterion (see \cite[Th.22.3]{Mat}), assertion
$\mathrm{(a)}_{n,j+1}$ will follow from :

\begin{claim}\label{cl_zumzum}
The induced morphism $\bar\bomega^*_n(p^jB/p^{j+1}B)\otimes_{W_{n,j}B}
W_{n,j}C\to\bar\bomega^*_n(p^jC/p^{j+1}C)$ is an isomorphism.
\end{claim}
\begin{pfclaim} Clearly
$\bar\bomega^*_n(p^jB/p^{j+1}B)$ is a $\bar\bomega^*_nB_0$-module
and $\bar\bomega^*_n(p^jC/p^{j+1}C)$ is a $C_0$-module, and notice
that $\bar\bomega^*_n(B_0)=\bar\bomega^*_0\circ\Phi^{n*}_{B_0}(B_0)$,
and likewise for $\bar\bomega^*_n(C_0)$. Then the stated morphism
agrees with the composition
$$
\begin{aligned}
\bar\bomega^*_n(p^jB/p^{j+1}B)\otimes_{W_{n,j}B}W_{n,j}C\isom&\,
\bar\bomega^*_0\Phi^{n*}_{B_0}(p^jB/p^{j+1}B)\otimes_{\bar\bomega^*_0B_0}
\bar\bomega^*_0B_0\otimes_{W_{n,j}B}W_{n,j}C \\
\isom&\,
\bar\bomega^*_0\Phi^{n*}_{B_0}(p^jB/p^{j+1}B)\otimes_{\bar\bomega^*_0B_0}
\bar\bomega^*_0C_0 \\
=&\,\bar\bomega_0^*(\Phi^{n*}_{B_0}(p^jB/p^{j+1}B)\otimes_{B_0}C_0) \\
\isom&\,\bar\bomega_0^*(\Phi^{n*}_{B_0}(p^jB/p^{j+1}B)
\otimes_{\Phi^{n*}_{B_0}B_0}\Phi^{n*}_{B_0}B_0\otimes_{B_0}C_0) \\
\isom&\,\bar\bomega_0^*(\Phi^{n*}_{B_0}(p^jB/p^{j+1}B)
\otimes_{\Phi^{n*}_{B_0}B_0}\Phi^{n*}_{B_0}C_0) \\
=&\,\bar\bomega_0^*\Phi^{n*}_{B_0}((p^jB/p^{j+1}B)\otimes_{B_0}C_0) \\
\isom&\,\bar\bomega_0^*\Phi^{n*}_{B_0}(p^jC/p^{j+1}C) \\
=&\,\bar\bomega_n^*(p^jC/p^{j+1}C)
\end{aligned}
$$
where the second isomorphism follows from $\mathrm{(b)}_{n,j}$,
the fourth follows from our assumption about $\Phi_{C_0/B_0}$,
and the fifth follows from the flatness of $f$.
\end{pfclaim}

Claim \ref{cl_zumzum} also implies that the diagram
$$
\xymatrix{ W_{n,j+1}B \ar[r] \ar[d]_{W_{n,j+1}f} &
W_{n,j}B \ar[d]^{W_{n,j}f} \\
W_{n,j+1}C \ar[r] & W_{n,j}C
}$$
is cocartesian; combining with $\mathrm{(b)}_{n,j}$, we see
that $\mathrm{(b)}_{n,j+1}$ holds. Thus, we conclude that
$\mathrm{(a)}_{n,j}$ and $\mathrm{(b)}_{n,j}$ hold for every
$j\in\N$; for large values of $j$, we deduce that
$\mathrm{(a)}_{n+1,0}$ and $\mathrm{(b)}_{n+1,0}$ both hold.
By induction on $n\in\N$, the assertion follows. This
concludes the proof of (i) in case $p^kA=0$. To prove
(ii) under the same assumption, we show first that the
natural morphism
\set\begin{equation}\label{eq_wiesia}
\bar\bomega^*_iB_0\otimes_{W_{n+1}B}W_{n+1}C\to\bar\bomega^*_iC_0
\end{equation}
is an isomorphism. Indeed, notice that
$\bar\bomega^*_iB_0=\bar\bomega^*_0(\Phi^{i*}_{B_0}B_0)$, and
likewise for $\bar\bomega^*_iC_0$. Then the foregoing morphism
is the composition :
$$
\begin{aligned}
\bar\bomega^*_iB_0\otimes_{W_{n+1}B}W_{n+1}C=&\,
\bar\bomega^*_0(\Phi^{i*}_{B_0}B_0)\otimes_{W_{n+1}B}W_{n+1}C \\
\isom&\,\bar\bomega^*_0(\Phi^{i*}_{B_0}B_0)\otimes_{\bar\bomega^*_0B_0}
\bar\bomega^*_0B_0\otimes_{W_{n+1}B}W_{n+1}C \\
\isom&\,\bar\bomega^*_0(\Phi^{i*}_{B_0}B_0)\otimes_{\bar\bomega^*_0B_0}
\bar\bomega^*_0C_0 \\
=&\,\bar\bomega^*_0(\Phi^{i*}_{B_0}B_0\otimes_{B_0}C_0) \\
\isom&\,\bar\bomega^*_0(\Phi^{i*}_{B_0}C_0) \\
=&\,\bar\bomega^*_iC_0
\end{aligned}
$$
where the second isomorphism follows from the foregoing
condition $\mathrm{(b)}_{n+1,0}$, and the third isomorphism
follows from our assumption on $\Phi_{C_0/B_0}$. Now, since
$p^kB=0$, the assertion follows from the isomorphism
\eqref{eq_wiesia} together with the following :

\begin{claim} Let $g:R\to S$ be a flat morphism of
$A^a$-algebras, $I\subset R$ a nilpotent ideal, and
suppose that $g\otimes_RR/I$ is an isomorphism. Then
$g$ is an isomorphism. 
\end{claim}
\begin{pfclaim} By a simple induction, we are reduced
to the case where $I^2=0$. Then we have the commutative
ladder with exact rows :
$$
\xymatrix{ 0 \ar[r] & I \ar[r] \ar[d] & R \ar[r] \ar[d] &
R/I \ar[r] \ar[d] & 0 \\
0 \ar[r] & I\otimes_RS \ar[r] & S \ar[r] & S/IS \ar[r] & 0
}$$
and a natural identification
$I\otimes_RS\isom I\otimes_{R/I}S/IS\isom I$. Thus, the
first and third vertical arrows are both isomorphisms,
and the same then follows for the middle one.
\end{pfclaim}

We consider next the case where $A$ is a general ring,
and $\fm$ fulfills condition $(\bB)$. We shall use the
following criterion :

\begin{claim}\label{cl_using-old-Ferrand}
Let $(V,\fm_V)$ be any basic setup, $R$ a
$(V,\fm_V)^a$-algebra and $M$ an $R$-module. Set
$$
M_\mathrm{tor}:=\bigcup_{n\in\N}\Ann_M(p^n)
\qquad
R_\mathrm{tor}:=\bigcup_{n\in\N}\Ann_R(p^n).
$$
Then $M$ is a flat $R$-module if and only if the
following three conditions
hold :
\begin{enumerate}
\alphaenu
\item
$M/p^iM$ is a flat $R/p^iR$-module for every $i\in\N$.
\item
$M\otimes_RR[p^{-1}]$ is a flat $R[p^{-1}]$-module.
\item
The natural morphism $R_\mathrm{tor}\otimes_RM\to M_\mathrm{tor}$
is an isomorphism.
\end{enumerate}
\end{claim}
\begin{pfclaim} Clearly if $M$ is flat, conditions (a)--(c)
hold. Thus, suppose that (a)--(c) hold; according to
\cite[Lemma 5.2.1]{Ga-Ra} it suffices to prove that
\set\begin{equation}\label{eq_from-Ferrand-I-guess}
\Tor^R_i(M,R/pR)=0
\qquad
\text{for $i=1,2$}.
\end{equation}
However, a standard calculation shows that
$$
\Tor^R_1(M,R/pR)=\frac{\Ann_M(p)}{\Ann_R(p)M}
\qquad\text{and}\qquad
\Tor^R_2(M,R/pR)=\Ker\,(\Ann_R(p)\otimes_RM\to M).
$$
Let us then consider the commutative ladder :
$$
\xymatrix{ 0 \ar[r] & \Ann_R(p)\otimes_RM \ar[r] \ar[d] &
R_\mathrm{tor}\otimes_RM \ar[r]^-p \ar[d] &
R_\mathrm{tor}\otimes_RM \ar[d] \\
0 \ar[r] & \Ann_M(p) \ar[r] & M_\mathrm{tor} \ar[r]^-p &
M_\mathrm{tor}
}$$
whose bottom horizontal row is exact, and whose central and
right vertical arrows are isomorphisms; in order to verify
\eqref{eq_from-Ferrand-I-guess}, it suffices therefore to
check that the top horizontal row is also exact; the latter
is the inductive limit of the system of morphisms :
$$
\Ann_R(p)\otimes_RM\to\Ann_R(p^n)\otimes_RM\to\Ann_R(p^n)\otimes_RM
\qquad
\text{for every $n\in\N$}.
$$
But clearly $\Ann_R(p^i)\otimes_RM=\Ann_R(p^i)\otimes_{R/p^nR}M/p^nM$
for every $i,n\in\N$ with $i\leq n$; in view of (a), the assertion
follows.
\end{pfclaim}

Now, example \ref{ex_localize-Witt}(ii), implies easily
that $W_{n+1}(C)[p^{-1}]$ is a flat $W_{n+1}(B)[p^{-1}]$-algebra.
Next, for every $k\in\N$ let $\pi_k:C\to C/p^kC$ be the
projection; from example \ref{ex_how-to-write-p-powers}(iii)
we get
$$
\Ker\,W_{n+1}(\pi_{n+k})\subset p^kW_{n+1}(C)
\qquad
\text{for every $k\in\N$}
$$
so that $\pi_{n+k}$ induces an isomorphism of $W_{n+1}A^a$-algebras
\set\begin{equation}\label{eq_double-quotients}
W_{n+1}(C)\otimes_\Z\Z/p^k\Z\isom
W_{n+1}(C/p^{n+k}C)\otimes_\Z\Z/p^k\Z
\qquad
\text{for every $k\in\N$}.
\end{equation}
But, by the foregoing case, we know already that
$W_{n+1}(C/p^{n+k}C)$ is a flat $W_{n+1}(B/p^{n+k}B)$-algebra
for every $k\in\N$, so we conclude that
$W_{n+1}(C)\otimes_\Z\Z/p^k\Z$ is a flat
$W_{n+1}(B)\otimes_\Z\Z/p^k\Z$-algebra for every $k\in\N$.
In light of claim \ref{cl_using-old-Ferrand}, we are
therefore reduced to checking that the natural morphism
$$
W_{n+1}(B)_\mathrm{tor}\otimes_{W_{n+1}B}W_{n+1}C\to W_{n+1}(C)_\mathrm{tor}
$$
is an isomorphism. More precisely, we shall show, by
descending induction on $i$, that :
\begin{itemize}
\item
the induced morphism
$(\bar V_iB)_\mathrm{tor}\otimes_{W_{n+1}B}W_{n+1}C\to
(\bar V_iC)_\mathrm{tor}$ is an isomorphism of $W_{n+1}B$-modules
for every $i=0,\dots,n+1$
\item
the induced morphism
$(\bar V_{i+1}B)_\mathrm{tor}\otimes_{W_{n+1}B}W_{n+1}C\to
(\bar V_iB)_\mathrm{tor}\otimes_{W_{n+1}B}W_{n+1}C$
is a monomorphism of $W_{n+1}B$-modules for $i=0,\dots,n$.
\end{itemize}
Indeed, both assertions are trivial for $i=n+1$. Suppose
that both assertions are already known for some strictly
positive integer $i\leq n+1$; from the commutative diagram
$$
\xymatrix{ (\bar V_iB)_\mathrm{tor}\otimes_{W_{n+1}B}W_{n+1}C
\ar[r] \ar[d] & (\bar V_iC)_\mathrm{tor} \ar[d] \\
(\bar V_{i-1}B)_\mathrm{tor}\otimes_{W_{n+1}B}W_{n+1}C
\ar[r] & (\bar V_{i-1}C)_\mathrm{tor} 
}$$
we then deduce that the left vertical arrow is a monomorphism,
{\em i.e.} the second assertion holds for $i-1$. Moreover,
in order to prove the first assertion for $i-1$, it suffices
to check that the induced morphism
$$
(\bar V_{i-1}B)_\mathrm{tor}/(\bar V_iB)_\mathrm{tor}
\otimes_{W_{n+1}B}W_{n+1}C\to
(\bar V_{i-1}C)_\mathrm{tor}/(\bar V_{i-1}C)_\mathrm{tor}
$$
is an isomorphism of $W_{n+1}B$-modules. On the other hand,
from example \ref{ex_localize-Witt}(i) we see that
$W_{n+1}(B)_\mathrm{tor}$ is the kernel of the natural
morphism $W_{n+1}B\to W_{n+1}(B[p^{-1}])$, and likewise
for $W_{n+1}(C)_\mathrm{tor}$; taking into account
\eqref{eq_from-cl-W_n-structure}, we deduce a natural
isomorphism :
$$
(\bar V_{i-1}B)_\mathrm{tor}/(\bar V_iB)_\mathrm{tor}\isom
(\bar\bomega_{i-1}^*B)_\mathrm{tor}
$$
of $W_{n+1}B$-modules, and likewise for
$(\bar V_{i-1}C)_\mathrm{tor}/(\bar V_{i-1}C)_\mathrm{tor}$.
We are then further reduced to checking that the induced
morphism
$$
\bar\bomega_{i-1}^*(\Ann_B(p^k))\otimes_{W_{n+1}(B)}W_{n+1}C\to
\bar\bomega_{i-1}^*\Ann_C(p^k)
$$
is an isomorphism for every $k\in\N$. However, the latter
is the composition of the isomorphisms
$$
\begin{aligned}
\bar\bomega_{i-1}^*(\Ann_B(p^k))\otimes_{W_{n+1}(B)}W_{n+1}C
\isom&\,\bar\bomega_{i-1}^*(\Ann_B(p^k))
\otimes_{W_{n+1}(B/p^{n+k}B)}W_{n+1}(C/p^{n+k}C) \\
\isom&\,\bar\bomega_{i-1}^*(\Ann_B(p^k))
\otimes_{\bar\bomega_{i-1}^*{(B/p^{n+k}B)}}\bar\bomega_{i-1}^*(C/p^{n+k}C) \\
=&\,\bar\bomega_{i-1}^*(\Ann_B(p^k)\otimes_BC) \\
\isom&\,\bar\bomega_{i-1}^*\Ann_C(p^k)
\end{aligned}
$$
where the first isomorphism is due to \eqref{eq_double-quotients},
the second one follows from the part of assertion (ii) that
we have already proved, and the third one follows from
the flatness of $f$.

To conclude the proof of (ii) for a general ring $A$,
we shall need the following :

\begin{claim}\label{cl_map-of-flat-mods}
Let $(V,\fm_V)$ be a basic setup, $R$ a $(V,\fm)^a$-algebra
and $\phi:M\to N$ a morphism of flat $R$-modules. Suppose that :
\begin{enumerate}
\alphaenu
\item
$\phi\otimes_\Z\F_p$ is an isomorphism of $R/pR$-modules.
\item
$\phi\otimes_\Z\Z[p^{-1}]$ is an isomorphism of
$R[p^{-1}]$-modules.
\end{enumerate}
Then $\phi$ is an isomorphism.
\end{claim}
\begin{pfclaim} Denote by $C_\bullet$ the complex of $R$-modules
$[M\xrightarrow{\phi}N]$, say with $M$ placed in degree $0$.
Consider as well the acyclic complex
$$
D_\bullet
\quad :\quad
0\to\Ann_R(p)\to R\xrightarrow{\ p\one_R\ }R\to R/pR\to 0.
$$
Since $C_\bullet$ is flat in every degree, the complex
$D_\bullet\otimes_RC_\bullet$ is still exact; the latter
is also the total complex of the double complex
$$
0\to E_\bullet:=\Ann_R(p)\otimes_{R/pR}(C_\bullet\otimes_\Z\F_p)
\to C_\bullet\xrightarrow{\ p\one_{C_\bullet}\ }C_\bullet\to
C_\bullet\otimes_\Z\F_p\to 0.
$$
On the other hand, assumption (a) says that
$C_\bullet\otimes_\Z\F_p$ is an acyclic complex of flat
$R/pR$-modules, hence also $E_\bullet$ is acyclic.
Summing up, we find that the morphism of complexes
$p\one_{C_\bullet}:C_\bullet\to C_\bullet$ is a quasi-isomorphism,
and hence induces isomorphisms
$p\cdot\one_{H_i(C_\bullet)}:H_i(C_\bullet)\isom H_i(C_\bullet)$
for $i=0,1$. In other words,
$H_i(C_\bullet)=H_i(C_\bullet)\otimes_\Z\Z[p^{-1}]$ for
$i=0,1$; combining with assumption (b), we conclude
that $\Ker\,\phi=\Coker\,\phi=0$, whence the claim.
\end{pfclaim}

Now, we need to show that the morphism
$$
\phi_f:\bar\bomega^*_iB\otimes_{W_{n+1}B}W_{n+1}C\to\bar\bomega^*_iC
$$
resulting from the diagram of (ii) is an isomorphism
of $W_{n+1}B$-modules; however, set as well
$g:=f\otimes_\Z\Z/p^{n+1}\Z:B/p^{n+1}B\to C/p^{n+1}C$;
according to the part of (ii) that has already been
proved, the corresponding morphism
$$
\phi_g:\bar\bomega^*_i(B/p^{n+1}B)
\otimes_{W_{n+1}(B/p^{n+1}B)}W_{n+1}(C/p^{n+1}C)\to
\bar\bomega^*_i(C/p^{n+1}C)
$$
is an isomorphism; on the other hand, in view of
\eqref{eq_double-quotients} the morphism $\phi_f\otimes_\Z\F_p$
is naturally identified with $\phi_g\otimes_\Z\F_p$, so
$\phi_f\otimes_\Z\F_p$ is an isomorphism as well. Next,
recall that the ghost map induces an isomorphism of
$W_{n+1}A^a$-algebras
$$
(W_{n+1}B)[p^{-1}]\isom B[p^{-1}]^{n+1}:=\prod_{i=0}^n
\bar\bomega^*_iB[p^{-1}]
$$
and likewise for $(W_{n+1}C)[p^{-1}]$ (example
\ref{ex_localize-Witt}(ii)); under these isomorphisms, the
morphism $\bar\bomega_i\otimes_\Z\Z[p^{-1}]$ on $(W_{n+1}B)[p^{-1}]$
is identified with the projection $\pi_i:B[p^{-1}]^{n+1}\to B[p^{-1}]$
on the $(i+1)$-th factor, and likewise for the corresponding
morphism on $(W_{n+1}C)[p^{-1}]$. Thus, $\phi_f\otimes_\Z\Z[p^{-1}]$
is identified with the natural isomorphism of $B[p^{-1}]^{n+1}$-modules
$$
\pi^*_iB[p^{-1}]\otimes_{B[p^{-1}]^{n+1}}C[p^{-1}]^{n+1}\isom
\pi^*_iC[p^{-1}].
$$
From claim \ref{cl_map-of-flat-mods} and (i) we deduce
that $\phi_f$ is an isomorphism, as required.
\end{proof}

\begin{corollary}\label{cor_Witt-and-etale-maps}
In the situation of theorem {\em\ref{th_Witt-and-etale-maps}},
for every $n\in\N$ we have :
\begin{enumerate}
\item
Every $B$-algebra $D$ induces a cocartesian diagram
of\/ $W_{n+1}A^a$-algebras :
$$
\xymatrix{ W_{n+1}B \ar[r] \ar[d] & W_{n+1}C \ar[d] \\
W_{n+1}D \ar[r] & W_{n+1}(C\otimes_BD).
}$$
\item
The $W_{n+1}B$-module $W_{n+1}C$ is almost finitely generated
(resp. almost finitely presented, resp. almost projective)
if and only if the same holds for the $B$-module $C$.
\item
Suppose that condition $(\bB)$ of \cite[\S2.1.6]{Ga-Ra} holds
for $\fm$. Then $W_{n+1}(f)$ is \'etale (resp. weakly \'etale)
if and only if the same holds for the morphism $f$.
\item
The projection $W_{n+2}C\to W_{n+1}C$ induces an isomorphism
of \/ $W_{n+1}B$-algebras :
$$
W_{n+2}C\otimes_{W_{n+2}B}W_{n+1}B\isom W_{n+1}C.
$$
\end{enumerate}
\end{corollary}
\begin{proof}(i): Notice that the natural morphism
$h:W_{n+1}C\otimes_{W_{n+1}B}W_{n+1}D\to W_{n+1}(C\otimes_BD)$
restricts to morphisms of $W_{n+1}B$-modules
$$
W_{n+1}C\otimes_{W_{n+1}B}\bar V_kD\to\bar V_k(C\otimes_BD)
\qquad
\text{for every $k=0,\dots,n$}.
$$
Taking into account \eqref{eq_from-cl-W_n-structure}, we
are then reduced to checking that $h$ induces
isomorphisms
$$
\gr_kh:W_{n+1}C\otimes_{W_{n+1}B}\bar\bomega^*_kD\to
\bar\bomega^*_k(C\otimes_BD)
\qquad
\text{for every $k=0,\dots,n$}.
$$
However, a simple inspection shows that $\gr_kh$
agrees with the composition
$$
W_{n+1}C\otimes_{W_{n+1}B}\bar\bomega^*_kD\isom
W_{n+1}C\otimes_{W_{n+1}B}\bar\bomega^*_kB\otimes_{\bar\bomega^*_kB}
\bar\bomega^*_kD\isom\bar\bomega^*_kC\otimes_{\bar\bomega^*_kB}
\bar\bomega^*_kD=\bar\bomega^*_k(C\otimes_BD)
$$
where the second isomorphism follows from theorem
\ref{th_Witt-and-etale-maps}(ii).

(iv): The proof of (i) also shows that the natural morphism
$$
W_{n+2}C\otimes_{W_{n+2}B}\bar V_{n+1}B\to\bar V_{n+1}C
$$
is an isomorphism. The assertion is an immediate consequence.

(ii): Let $\psi:W_{n+1}B\to B^{n+1}:=\prod_{i=0}^n\bar\bomega^*_iB$
be the morphism of $W_{n+1}A^a$-algebras whose composition with
the projection $B^{n+1}\to\bar\bomega^*_iB$ is the morphism
$\bar\bomega_i$, for $i=0,\dots,n$; denote by $R\subset B^{n+1}$
the image of $\psi$. Let also $\bP$ be one of the properties :
``almost finitely generated'', ``almost finitely presented'',
or ``almost projective'', and suppose that $f$ enjoys $\bP$;
in view of theorem \ref{th_Witt-and-etale-maps}(i), proposition
\ref{prop_int-with-nilker}(i) and \cite[Lemma 3.2.25]{Ga-Ra}, it
suffices to prove that the $R$-module $W_{n+1}C\otimes_{W_{n+1}B}R$
enjoys property $\bP$. Next, by proposition
\ref{prop_int-with-nilker}(ii) we know that
$$
p^{n+1}B^{n+1}\subset R
\qquad\text{and we let}\qquad
R_0:=R/p^{n+1}B^{n+1}
\qquad
B_0^{n+1}:=B^{n+1}/pB^{n+1}.
$$
With this notation, it is easily seen that there exists
a well defined morphism of $W_{n+1}B$-algebras
$$
\phi:B_0\to B_0^{n+1}
\qquad\text{such that}\qquad
\pi_i\circ\phi=\Phi^i_{B_0}:B_0\to
\bar\bomega_i^*B_0=\Phi^{i*}_{B_0}B_0
\qquad
\text{for $i=0,\dots,n$}
$$
where $\pi_i:B_0^{n+1}\to\bar\bomega_i^*B_0$ is the projection
on the $(i+1)$-th factor. Arguing as in the proof of proposition
\ref{prop_alm-str-on-Witt}, we get cartesian diagrams
$$
\xymatrix{ S \ar[r] \ar[d] & B^{n+1} \ar[d] &
R \ar[r] \ar[d] & R_0 \ar[d] \\
B_0 \ar[r]^-\phi & B_0^{n+1}  & S \ar[r] & S_0:=S/p^{n+1}B^{n+1}.
}$$
In case $\bP$ is ``almost finitely generated'' or ``almost
finitely presented'' we may then apply
\cite[Rem.3.2.26(i) and Lemma 3.4.18(i)]{Ga-Ra}
to the second of these cartesian diagrams, and thereby
reduce to checking property $\bP$ for the $(S\times R_0)$-module
$W_{n+1}C\otimes_{W_{n+1}B}(S\times R_0)$. In case $\bP$ is the
property ``almost projective'', we may argue as in the proof
of \cite[Prop.3.4.21]{Ga-Ra} to achieve the same reduction.
But notice as well that the restriction
$R_0\to\bar\bomega^*_0B_0$ of $\pi_0$ is an epimorphism (on the
underlying modules), and has nilpotent kernel (see the proof
of proposition \ref{prop_alm-str-on-Witt}), so the $R_0$-module
$W_{n+1}C\otimes_{W_{n+1}B}R_0$ enjoys $\bP$ if and only the same
holds for the $\bar\bomega^*_0B_0$-module
$W_{n+1}C\otimes_{W_{n+1}B}\bar\bomega^*_0B_0$. Next, by invoking
again \cite[Rem.3.2.26(i) and Lemma 3.4.18(i)]{Ga-Ra} or
the proof of \cite[Prop.3.4.21]{Ga-Ra} to the first of the
two above diagrams we deduce that the $S$-module
$W_{n+1}C\otimes_{W_{n+1}B}S$ enjoys $\bP$ if and only if
the same holds for the $(\bar\bomega^*_0B_0\times B^{n+1})$-module
$W_{n+1}C\otimes_{W_{n+1}B}(\bar\bomega^*_0B_0\times B^{n+1})$. Since
$\bar\bomega^*_0B_0$ is a quotient of $B^{n+1}$, we are
finally reduced to checking that the $B^{n+1}$-module
$W_{n+1}C\otimes_{W_{n+1}B}B^{n+1}$ enjoys $\bP$; but according
to theorem \ref{th_Witt-and-etale-maps}(ii), the latter is
isomorphic to the $B^{n+1}$-module $\prod_{i=0}^n\bar\bomega_i^*C$,
so the assertion is clear.

(iii): Suppose first that $f$ is weakly \'etale, so that the
multiplication morphism $\mu_{C/B}:C\otimes_BC\to C$ is flat;
since it is also an epimorphism on the underlying $B$-modules,
it follows easily that $\mu_{C/B}$ is weakly \'etale. Especially,
$\Phi_{C/C\otimes_BC}$ is an isomorphism, by remark
\ref{rem_explain-Frob} and our assumption on $\fm$.
Then theorem \ref{th_Witt-and-etale-maps} says that
$W_{n+1}(\mu_{C/C\otimes_BC})$ is flat. Lastly, according to
(i), the multiplication morphism $\mu_{W_{n+1}C/W_{n+1}B}$ of
$W_{n+1}C$ factors through $W_{n+1}(\mu_{C/C\otimes_BC})$ and an
isomorphism
\set\begin{equation}\label{eq_Witt-diagonal}
W_{n+1}C\otimes_{W_{n+1}B}W_{n+1}C\isom W_{n+1}(C\otimes_BC).
\end{equation}
Hence $\mu_{W_{n+1}C/W_{n+1}B}$ is flat, and combining with
theorem \ref{th_Witt-and-etale-maps}, we conclude that
$W_{n+1}f$ is weakly \'etale. In case $f$ is \'etale, we
deduce from (ii), the isomorphism \eqref{eq_Witt-diagonal},
and \cite[Prop.2.4.18]{Ga-Ra} that $W_{n+1}f$ is \'etale
as well. Conversely, (iv) implies that if $W_{n+1}f$ is
weakly \'etale (resp. \'etale), then the same holds for $f$.
\end{proof}

\begin{remark}
(i)\ \
It is shown in \cite[1.5.8]{Il2} that if $f:B\to C$ is a
homomorphism of $\F_p$-algebras and $n\in\N$ any integer,
then $W_{n+1}f$ is \'etale (in the usual sense of \cite{EGA4},
which includes the condition that $f$ is finitely presented),
if and only if the same holds for $f$. This result can be
deduced from corollary \ref{cor_Witt-and-etale-maps}(iii)
and extended to arbitrary rings, as follows. Resume the notation
of the proof of the corollary, and suppose that $\fm=A$, so we
are dealing with the ``classical limit'' for which almost rings
are just usual rings; then, $W_{n+1}f$ is a weakly \'etale ring
homomorphism, and by \cite[\S3.4.44]{Ga-Ra} we know already
that $W_{n+1}f$ is finitely presented if and only if the same
holds for $W_{n+1}f\otimes_{W_{n+1}B}R$. We claim next that
$W_{n+1}f\otimes_{W_{n+1}B}R$ is finitely presented if and
only if the same holds for $W_{n+1}f\otimes_{W_{n+1}B}(R_0\times S)$.
For the proof, arguing as in \cite[\S3.4.44]{Ga-Ra} we
reduce to checking that the natural morphism
$R\to R_0\times S$ is of universal effective descent
for the fibred category of weakly \'etale morphisms of
rings. The latter assertion is already known by corollary
\ref{cor_instead-of-3.2.1}. By the same token,
$W_{n+1}f\otimes_{W_{n+1}B}S$ is finitely presented if and
only if the same holds for
$W_{n+1}f\otimes_{W_{n+1}B}(B^{n+1}\times B_0)$. Moreover, the
surjective ring homomorphism $R_0\to B_0$ has nilpotent
kernel, so again by \cite[\S3.4.44]{Ga-Ra} we know that
$W_{n+1}f\otimes_{W_{n+1}B}R_0$ is finitely presented
if and only if the same holds for $W_{n+1}f\otimes_{W_{n+1}B}B_0$.
Since $B_0$ is a quotient of $B^{n+1}$, we are finally reduced
to checking that $W_{n+1}f\otimes_{W_{n+1}B}B^{n+1}$ is finitely
presented; but it was observed in the proof of corollary
\ref{cor_Witt-and-etale-maps}(ii) that the latter is isomorphic
to the $B^{n+1}$-module $\prod_{i=0}^n\bar\bomega_i^*C$, so the
assertion is clear.

(ii)\ \
A version of corollary \ref{cor_Witt-and-etale-maps}(iii)
for truncated big Witt vectors of usual rings is found
in \cite[Th.2.4]{vdKal}. Another proof and generalization
is given by \cite[Th.9.2]{Borg}.
\end{remark}

\subsection{Complements : locally measurable algebras}
\label{sec_loc-measur-algebras}
This section studies the global counterpart of the class of
measurable algebras introduced in section \ref{sec_norm-lengths}.
To begin with -- and until \eqref{subsec_end-of-digres} -- we
consider an arbitrary valued field $(K,|\cdot|)$, and we resume
the notation of \eqref{sec_sch-val-rings} and \eqref{sec_norm-lengths}.
Our first result is the following generalization of proposition
\ref{prop_Gruson-out} :

\begin{proposition}\label{prop_Gruson-further}
Let $A$ be a measurable $K^+$-algebra, $M$ a $K^+$-flat and
finitely generated $A$-module. Then $M$ is a finitely
presented $A$-module.
\end{proposition}
\begin{proof} Let $\Sigma$ be a finite system of generators
for $M$; also let us write $A$ as the colimit of a filtered
system $(A_i~|~i\in I)$ of finitely presented $K^+$-algebras,
with \'etale transition maps. For every $i\in I$, let $M_i$
be the $A_i$-submodule of $M$ generated by $\Sigma$; notice
that $M_i$ is still $K^+$-flat, hence it is a finitely presented
$A_i$-module (proposition \ref{prop_Gruson-out}). We may then
write $M$ as the colimit of the filtered system
$(M_i\otimes_{A_i}A~|~i\in I)$ of finitely presented $A$-modules,
with surjective transition maps. Moreover, since $\bar A:=A/\fm_KA$
is noetherian (lemma \ref{lem_lucilla}(i)), there exists $i\in I$
such that $M_j\otimes_{A_j}\bar A=M/\fm_KM$ for every $j\geq i$.
Since the ring homomorphism $\bar A_j:=A_j/\fm_KA_j\to\bar A$
is faithfully flat, we deduce that
\set\begin{equation}\label{eq_deduce-that}
M_i\otimes_{A_i}\bar A_j=M_j\otimes_{A_j}\bar A_j
\qquad
\text{for every $j\geq i$}.
\end{equation}
Consider the short exact sequence
$$
C\quad :\quad
0\to N\to M_i\otimes_{A_i}A_j\to M_j\to 0.
$$
Since $M_i$ is $K^+$-flat, the same holds for $N$, and the
latter is a finitely generated $A_j$-module, since both $M_j$
and $M_i\otimes_{A_i}A_j$ are finitely presented; therefore,
the sequence $C\otimes_{K^+}\kappa$ is still exact. Taking
into account \eqref{eq_deduce-that}, we see that $N/\fm_KN=0$.
By Nakayama's lemma, it follows that $N=0$, and finally,
$M=M_i\otimes_{A_i}A$ is finitely presented, as stated.
\end{proof}

\begin{definition} Let $A$ be a $K^+$-algebra. Set $X:=\Spec\,A$,
$S:=\Spec\,K^+$, and denote by $f:X\to S$ the structure morphism.
We say that $A$ is {\em locally measurable}, if the following
holds. For every $x\in X$ and every point $\xi$ of $X$ localized
at $x$, the strict henselization of $A$ at $\xi$ is a measurable
$\cO_{\!S,f(x)}$-algebra.
\end{definition}

\begin{remark}\label{rem_local-measura}
Let $A$ be a local $K^+$-algebra, $A^\sh$ the strict
henselization of $A$ at a geometric point localized at the
closed point, and $M$ an $A$-module.

(i)\ \
Clearly, $A$ is locally measurable if and only if $A^\sh$
is measurable.

(ii)\ \
However, if $A^\sh$ is measurable, it does not necessarily
follow that $A$ is measurable.

(iii)\ \
On the other hand, if $A$ is a normal local domain, one can
show that $A$ is measurable if and only if the same holds
for $A^\sh$.

(iv)\ \
Suppose that $A$ is local and locally measurable. Since the
natural map $A\to A^\sh$ is faithfully flat, and since every
measurable $K^+$-algebra is a coherent ring, it is easily seen
that $A$ is coherent. If furthermore, the structure map
$K^+\to A$ is local, lemma \ref{lem_lucilla}(i) implies that
$A/\fm_KA$ is a noetherian ring.
\end{remark}

\begin{remark}
Let $A$ be a locally measurable $K^+$-algebra. Then, for
every finitely generated ideal $I\subset A$, the $K^+$-algebra
$A/I$ is also locally measurable. Indeed, for every geometric
point $\xi$ of $\Spec\,A$, let $A^\sh_\xi$ (resp. $(A/I)^\sh_\xi$)
be the strict henselization of $A$ (resp. of $A/I$) at $\xi$.
Then the natural map $A^\sh_\xi/IA^\sh_\xi\to(A/I)_\xi^\sh$
is an isomorphism for every such $\xi$
(\cite[Ch.IV, Prop.18.8.10]{EGA4}), so the assertion follows
from lemma \ref{lem_lucilla}(iv).
\end{remark}

\begin{definition}\label{def_breaks}
Let $A$ be a $K^+$-algebra, $M$ an $A$-module, and
$\gamma\in\log\Gamma^+$ any element.
\begin{enumerate}
\item
A {\em $K^+$-flattening sequence\/} for the $A$-module $M$ is
a finite sequence $\underline b:=(b_0,\dots,b_n)$ of elements
of $K^+$ such that
\begin{enumerate}
\item
$\log|b_{i+1}|>\log|b_i|$ for every $i=0,\dots,n-1$, $b_0=1$
and $b_n=0$.
\item
$b_iM/b_{i+1}M$ is a $K^+/b_i^{-1}b_{i+1}K^+$-flat module,
for every $i=0,\dots,n-1$.
\end{enumerate}
We say that a $K^+$-flattening sequence $\underline b$ for $M$
is {\em minimal}, if no proper subsequence of $\underline b$ is
still flattening for $M$.
\item
Say that $\gamma=\log|b|$ for some $b\in K^+$, and let
$b_M:M\to bM$ be the map given by the rule : $m\mapsto bm$,
for every $m\in M$. We say that $\gamma$ {\em breaks\/} $M$
if the $\kappa$-linear map
$$
b_M\otimes_{K^+}\kappa:
M\otimes_{K^+}\kappa\to bM\otimes_{K^+}\kappa
$$
is not an isomorphism.
\end{enumerate}
\end{definition}

\begin{remark}\label{rem_breaks}
Let $A$ be a $K^+$-algebra, $M$ an $A$-module.

(i)\ \
Suppose that $\gamma\in\log\Gamma^+$ breaks $M$, and say that
$\gamma=\log|c|$ for some $c\in K^+$. Clearly, for every
$b\in cK^+$, the map $b_M$ factors through $c_M$. We deduce
that every $\gamma'\in\log\Gamma^+$ with $\gamma'\geq\gamma$
also breaks $M$.

(ii)\ \
For given $b\in K^+$, suppose that $M$ is a flat
$K^+/bK^+$-module. Then we claim that no $\gamma<\log|b|$
breaks $M$. Indeed, for such $\gamma$ pick $c\in K^+$ with
$\log|c|=\gamma$, and set $W:=K^+/bK^+$; notice that the map
$c_W:W\to cW$ induces an isomorphism
$c_W\otimes_W\one_\kappa:\kappa\isom cW\otimes_W\kappa$
(notation of definition \ref{def_breaks}(ii)), whence an
isomorphism 
$$
\one_M\otimes_Wc_W\otimes_W\one_\kappa:
M\otimes_W\kappa\isom M\otimes_W(cW\otimes_W\kappa).
$$
Since $M$ is a flat $W$-module, the natural map
$M\otimes_W cW\to cM$ is an isomorphism, and the
resulting isomorphism $M\otimes_W\kappa\isom cM\otimes_W\kappa$
is naturally identified with $c_M$, whence the contention.

(iii)\ \
Let $\underline b:=(b_0,\dots,b_n)$ be a sequence of elements
of $K^+$ fulfilling condition (a) of definition
\ref{def_breaks}(i), and suppose that $\underline b$ admits
a subsequence that is $K^+$-flattening for $M$. Then
$\underline b$ is $K^+$-flattening for $M$ as well. Indeed,
an easy induction reduces the contention to the following.
Let $b,c\in K^+$ be any two elements such that $\log|c|<\log|b|$,
and suppose that $M$ is $K^+/bK^+$-flat; then $M/cM$ is
$K^+/cK^+$-flat, and $cM$ is $K^+/c^{-1}bK^+$-flat. Of this
two assertion, the first is trivial; to show the second,
recall that $M$ can be written as the colimit of a filtered
system of free $K^+/bK^+$-modules (\cite[Ch.I, Th.1.2]{La}).
Then we may assume that $M$ is free, in which case the
assertion is easily verified.
\end{remark}

\begin{lemma}\label{lem_flattening-strata}
Let $A$ be a $K^+$-algebra, $M$ a coherent $A$-module that
admits a $K^+$-flattening sequence, and suppose that the
Jacobson radical of $A$ contains $\fm_KA$. Then we have :
\begin{enumerate}
\item
$M$ admits a minimal $K^+$-flattening sequence, unique up to
units of $K^+$.
\item
Let $(b_0,\dots,b_n)$ be a minimal $K^+$-flattening sequence for
$M$. For every $\gamma\in\log\Gamma^+$ and every $i=0,\dots,n-1$,
the following conditions are equivalent :
\begin{enumerate}
\item
$\gamma$ breaks $b_iM$.
\item
$\gamma\geq\log|b_i^{-1}b_{i+1}|$. (As usual, we set
$\log|0|:=+\infty$ : see \eqref{sec_norm-lengths}.)
\end{enumerate}
\end{enumerate}
\end{lemma}
\begin{proof} It is clear that $M$ admits a minimal
$K^+$-flattening sequence $(1,b_1,\dots,b_{n-1},0)$, and (i)
asserts that every other such minimal sequence is of the type
$(1,u_1b_1,\dots,u_{n-1}b_{n-1},0)$ for some elements
$u_1,\dots,u_{n-1}\in(K^+)^\times$. However, notice that
the condition of (ii) characterize uniquely such a sequence,
so we only have to show that (ii) holds. We may also assume
that $n>1$; indeed, when $n=1$, the module $M$ is $K^+$-flat,
and the assertion is immediate. Moreover, since $M$ is coherent,
the same holds for $bM$, for every $b\in K^+$. Thus, an easy
induction reduces to showing the equivalence of conditions
(ii.a) and (ii.b) for $i=0$.

However, notice that if $\gamma$ breaks $M$, then it obviously
also breaks $M/b_1M$; in light of  remark \ref{rem_breaks}(ii),
it follows already that (a)$\Rightarrow$(b). Therefore, taking into
account remark \ref{rem_breaks}(i), it remains only to show
that $\log|b_1|$ breaks $M$.

If $n=2$, then $b_1M$ is a flat $K^+$-module; let
$N:=\Ker\,b_{1,M}$. We have already observed that $b_1M$ is a
finitely presented $A$-module, hence $N$ is a finitely generated
$A$-module, and $N\otimes_{K^+}\kappa$ is the kernel of
$b_{1,M}\otimes_{K^+}\kappa$. Suppose that $\log|b_1|$ does not
break $M$; then $N\otimes_{K^+}\kappa$ vanishes, and then $N=0$,
by Nakayama's lemma. It follows that $M$ is $K^+$-flat, which
contradicts the minimality of the sequence $(1,b_1,0)$.

Next, we consider the case where $n>2$ and set $M':=M/b_2M$,
$b:=b^{-1}_1b_2$; since $b_2 M\subset b_1\fm M$, it suffices
to show that $\log|b_1|$ breaks $M'$, hence we may assume that
$b_2M=0$, $b_1M$ is a flat $K^+/bK^+$-module, and the flattening
sequence is $(1,b_1,b_2,0)$. We remark :

\begin{claim}\label{cl_rapido}
If $\log|b_1|$ does not break $M$, the map
$b_{1,M/bM}:M/bM\to b_1M$ is an isomorphism of $K^+/bK^+$-modules.
\end{claim}
\begin{pfclaim} Indeed, let $N:=\Ker\,b_{1,M/bM}$; since $b_1M$
is $K^+/bK^+$-flat, $N\otimes_{K^+}\kappa$ is the kernel of
$b_{1,M/bM}\otimes_{K^+}\kappa$, and the latter vanishes if
$\log|b_1|$ does not break $M$; on the other hand, $N$ is a
finitely generated $A$-module, hence $N=0$, by Nakayama's lemma.
\end{pfclaim}

We shall use the following variant of the local flatness
criterion :

\begin{claim}\label{cl_Ferrando}
Let $R$ be a ring, $I_1,I_2\subset A$ two ideals, and $M$ an
$R$-module such that :
\begin{enumerate}
\alphaenu
\item
$I_1I_2=0$
\item
$M/I_iM$ is $R/I_i$-flat for $i=1,2$
\item
the natural map $I_1\otimes_RM/I_2M\to I_1M$ is an isomorphism.
\end{enumerate}
Then $M$ is a flat $R$-module.
\end{claim}
\begin{pfclaim} Set $I_3:=I_1\cap I_2$; to begin with,
\cite[Lemma 3.4.18]{Ga-Ra} and (b) imply that $M/I_3M$ is a flat
$R/I_3$-module. We have an obviously commutative diagram
$$
\xymatrix{
I_3\otimes_RM/I_3M \ar[rr]^-\alpha \ar[d]_\beta & & I_3M \ar[d]^\gamma \\
I_3\otimes_{R/I_2}M/I_2M \ar[r]^-\delta &
I_1\otimes_{R/I_2}M/I_2M \ar[r]^-\tau & I_1M
}$$
of $R$-linear maps. From (a) we see that $I_3^2=0$, and it
follows easily that $\beta$ is an isomorphism; $\gamma$ is
clearly injective, and the same holds for $\delta$, since
$M/I_2M$ is a flat $R/I_2$-module; by the same token, $\tau$
is an isomorphism. We conclude that $\alpha$ is an isomorphism,
and then the local flatness criterion (\cite[Th.22.3]{Mat})
yields the contention.
\end{pfclaim}

We shall apply claim \ref{cl_Ferrando} with $R:=K^+/b_2K^+$,
$I_1:=b_1R$, $I_2:=bR$. Indeed, condition (a) obviously holds
with these choices; by assumption, $M/I_1M$ is a $R/I_1$-flat
module, and if $\log|b_1|$ does not break $M$, claim
\ref{cl_rapido} implies that $M/I_2M$ is a flat $R/I_2$-module.
Lastly, condition (c) is equivalent to claim \ref{cl_rapido}.
Summing up, we have shown that if $\log|b_1|$ does not break
$M$, then $M$ is a flat $K^+/b_2$-module, contradicting again
the minimality of our flattening sequence.
\end{proof}

\begin{proposition}\label{prop_flattening-strata}
Let $A$ be a $K^+$-algebra, $M$ a finitely presented $A$-module,
and suppose that either one of the following two conditions holds :
\begin{enumerate}
\alphaenu
\item
$A$ is an essentially finitely presented $K^+$-algebra.
\item
$A$ is a local and locally measurable $K^+$-algebra.
\end{enumerate}
Then $M$ admits a $K^+$-flattening sequence.
\end{proposition}
\begin{proof} Suppose first that (b) holds, and let $A^\sh$ be
the strict henselization of $A$ at a geometric point localized
at the closed point. Since $A^\sh$ is a faithfully flat $A$-algebra,
it suffices to show that $M\otimes_AA^\sh$ admits a $K^+$-flattening
sequence. Hence, we may assume that $A$ is measurable, in which case
we may find a finitely presented $K^+$-algebra $A_0$ with an
ind-\'etale map $A_0\to A$ and a finitely presented $A_0$-module
$M_0$ with an isomorphism $M_0\otimes_{A_0}A\isom M$ of $A$-modules.
We may then replace $A$ by $A_0$, $M$ by $M_0$, and therefore
assume that $A$ is finitely presented over $K^+$, especially,
we are reduced to showing the assertion in the case where (a) holds.
To this aim, let us remark :

\begin{claim}\label{cl_capeesh}
Let $B:=\bigoplus_{n\in\N}B_n$ be a $\N$-graded finitely presented
$K^+$-algebra with $B_0=K^+$, and $N:=\bigoplus_{n\in\N}N_n$ a
$\N$-graded finitely presented $B$-module. We have :
\begin{enumerate}
\item
The $K^+$-module $N_n$ is finitely presented, for every $n\in\N$.
\item
For every $n\in\N$, let $(\gamma_{n,i}~|~i\in\N)$ be the sequence
of elementary divisors of $N_n$ (see \eqref{subsec_element-divs}).
Then $\Gamma(N):=\{\gamma_{n,i}~|~n,i\in\N\}$ is a finite set.
\item
$N$ admits a $K^+$-flattening sequence.
\end{enumerate}
\end{claim}
\begin{pfclaim}(i) is just a special case of proposition
\ref{prop_four-year-later}(iii).

(ii): Let $\cQ$ be the set of all finitely presented graded
quotients $Q$ of the $B$-modules $N$, for which $\Gamma(Q)$
is infinite. We have to show that $\cQ=\emptyset$. However,
for every $Q\in\cQ$, the $B\otimes_{K^+}\kappa$-module
$\bar Q:=Q\otimes_{K^+}\kappa$ is a quotient of
$\bar N:=N\otimes_{K^+}\kappa$; since $B\otimes_{K^+}\kappa$
is a noetherian ring, the set $\bar\cQ:=\{\bar Q~|~Q\in\cQ\}$
admits minimal elements, if it is not empty. In the latter
case, we may then replace $N$ by any $Q_0\in\cQ$ such that
$\bar Q_0$ is a minimal element of $\bar\cQ$, and assume that,
for every graded quotient $Q$ of $N$, either $\bar Q=\bar N$,
or else $\Gamma(Q)$ is a finite set. Now, for every $n\in\N$,
let $\gamma_n$ be the minimal non-zero elementary divisor of
$N_n$, and pick $a_n\in\fm_K$ with $\log|a_n|=\gamma_n$ (if
$N_n=0$, set $a_n=0$). Let $I\subset K^+$ be the ideal generated
by $\{a_n~|~n\in\N\}$; it is easily seen that $N_n/IN_n$ is
a free $K^+/I$-module for every $n\in\N$, hence $N/IN$ is
a $K^+/I$-flat finitely presented $B/IB$-module. We may then
find $a\in I$ such that $N/aN$ is already a $K^+/aK^+$-flat
$B/aB$-module (\cite[Ch.IV, Cor.11.2.6.1]{EGAIV-3}), so
$N_n/aN_n$ is a free $K^+/aK^+$-module for every $n\in\N$,
and we easily deduce that $\log|a|\leq\gamma_n$ for every
$n\in\N$, {\em i.e.} $I=aK^+$. Notice that $aN$ is a finitely
presented $B$-module (corollary \ref{cor_coherence}); it
follows that
$$
\Gamma(N)=\{\gamma+\log|a|~|~\gamma\in\Gamma(aN)\}\cup\{0\}.
$$
Especially, $\Gamma(aN)$ is an infinite set. On the other
hand, there exists some $n\in\N$ such that
$\dim_\kappa(aN_n)\otimes_{K^+}\kappa<
\dim_\kappa N_n\otimes_{K^+}\kappa$, so
$(aN)\otimes_{K^+}\kappa$ is a proper quotient of $\bar N$,
a contradiction.

(iii): Let $|b_0|,\dots,|b_{n-1}|$ be the finitely many elements
of $\Gamma(N)$ (for suitable $b_1,\dots,b_n\in K^+$), and set
$b_n:=0$; after permutation, we may assume that
$b_{i+1}\in b_i\fm_K$ for every $i=0,\dots,n-1$, and $b_0=1$.
Then we claim that $(b_0,\dots,b_n)$ is a $K^+$-flattening
sequence for $N$, {\em i.e.} $b_iN_k/b_{i+1}N_k$ is
$K^+/b_i^{-1}b_{i+1}K^+$-flat for every $i=0,\dots,n$ and every
$k\in\N$. However, say that $(\log|c_j|~|~j\in\N)$ is the
sequence of elementary divisors of $N_k$; we are reduced
to checking that $(c_jK^++b_iK^+)/(c_jK^++b_{i+1}K^+)$ is a flat
$K^+/b_i^{-1}b_{i+1}K^+$-module for every $j\in\N$. However,
by construction we have either $c_jK^+\subset b_{i+1}K^+$, or
$b_iK^+\subset c_jK^+$; in either case the assertion is clear.
\end{pfclaim}

Let now $A$ be an arbitrary essentially finitely presented
$K^+$-algebra, and $M$ an arbitrary finitely presented $A$-module.
We easily reduce to the case where $A=K^+[T_1,,\dots,T_r]$ is a
free polynomial $K^+$-algebra. In this case, we define a filtration
$\Fil_\bullet A$ on $A$, by declaring that $\Fil_k A$ is
the $K^+$-submodule of all polynomials of total degree $\leq k$,
for every $k\in\N$; then $R:=\sR(A,\Fil_\bullet A)_\bullet$
is a free polynomial $K^+$-algebra as well (see example
\ref{ex_Rees-free}). Let
$$
L_1\xrightarrow{\ \phi\ }L_0\to M
$$
be a presentation of $M$ as quotient of free $A$-modules
of finite rank. Let
$$
\be:=(e_1,\dots,e_n)
\qquad
\text{(resp.\ $\bff:=(f_1,\dots,f_m)$)}
$$
be a basis of $L_0$ (resp. of $L_1$); we endow $L_0$ with
the good $(A,\Fil_\bullet A)$-filtration $\Fil_\bullet L_0$
associated with the pair $(\be,(1,\dots,1))$ as in
\eqref{subsec_good-filtr} (this means that $e_i\in\Fil_1L_0$
for every $i=1,\dots,n$). Also, for every $i=1,\dots,m$,
pick $j_i\in\N$ such that $\phi(f_i)\in\Fil_{j_i}L_0$, and
endow $L_1$ with the good $(A,\Fil_\bullet A)$-filtration
$\Fil_\bullet L_1$ associated with the pair
$(\bff,(j_1,\dots,j_m))$. Set $L'_i:=\sR(L_i,\Fil_\bullet L_i)$,
and notice that $L'_i$ is a free $R$-module of finite rank,
for $i=0,1$. With these choices, $\phi$ is a map of filtered
$A$-modules, and there follows an $R$-linear map of
$\N$-graded $R$-modules $\sR(\phi)_\bullet:L'_1\to L'_0$,
whose cokernel is a $\N$-graded finitely presented
$R$-module $N_\bullet$. By inspecting the construction,
it is easily seen that the inclusion
$\Fil_kL_0\subset\Fil_{k+1}L_0$ induces a $K^+$-linear map
$N_k\to N_{k+1}$, for every $k\in\N$, as well as an
isomorphism of $K^+$-modules
$$
\colim_{k\in\N}N_k\isom M.
$$
On the other hand, claim \ref{cl_capeesh} ensures that $N_\bullet$
admits a $K^+$-flattening sequence. We easily deduce that the same
sequence is also flattening for $M$.
\end{proof}

\begin{corollary}\label{cor_flatness-mod-I}
Let $A$ be a local and locally measurable $K^+$-algebra,
$I\subset K^+$ any ideal, $M$ a finitely generated $A/IA$-module,
and suppose that the structure map $K^+\to A$ is local. Then the
following conditions are equivalent :
\begin{enumerate}
\alphaenu
\item
$M$ is a flat $K^+/I$-module.
\item
For every $c\in K^+$ such that $I\subset c\cdot\fm_K$, the value
$\log|c|$ does not break $M$.
\item
$M$ is a $K^+/I$-flat finitely presented $A/IA$-module.
\end{enumerate}
\end{corollary}
\begin{proof} (a)$\Rightarrow$(b): in light of remark
\ref{rem_breaks}(ii), it suffices to remark that
$M\otimes_{K^+}/bK^+$ is a flat $K^+/bK^+$-module, for
every $b\in K^+$ such that $I\subset K^+b$.

(b)$\Rightarrow$(c): Let us write $M$ as the colimit of a
filtered system $(M_\lambda~|~\lambda\in\Lambda)$ of finitely
presented $A$-modules, with surjective transition maps.
Then each $M_\lambda$ is a coherent $A$-module, and $A/\fm_KA$
is noetherian (remark \ref{rem_local-measura}(iv)), so we may
also assume that the induced maps
$\phi_\lambda:M_\lambda/\fm_KM_\lambda\to M/\fm_KM$ are
isomorphisms for every $\lambda\in\Lambda$. In view of the
commutative diagram
$$
\xymatrix{
M_\lambda\otimes_{K^+}\kappa
\ar[r]^-{c_{M_\lambda}} \ar[d]_{\phi_\lambda} &
cM_\lambda\otimes_{K^+}\kappa \ar[d] \\
M\otimes_{K^+}\kappa \ar[r]^-{c_M} & cM\otimes_{K^+}\kappa
}$$
it easily follows that $\log|c|$ does not break $M_\lambda$,
for any $\lambda\in\Lambda$ and any $c\in K^+$ such that
$I\subset c\cdot\fm_K$. Thus, if $(b_0,\dots,b_n)$ is the
minimal flattening sequence for $M_\lambda$, we see that
$b_1\in I$ (lemma \ref{lem_flattening-strata}(ii)), and
therefore $M_\lambda/IM_\lambda$ is a flat $K^+/I$-module
for every $\lambda\in\Lambda$.
Now, for $\lambda,\mu\in\Lambda$ with $\mu\geq\lambda$, let
$\psi_{\lambda\mu}:M_\lambda\to M_\mu$ be the transition map,
and set
$N_{\lambda\mu}:=\Ker(\psi_{\lambda\mu}\otimes_{K^+}K^+/I)$;
it follows that $N_{\lambda\mu}\otimes_{K^+}\kappa$ is the
kernel of $\psi_{\lambda\mu}\otimes_{K^+}\kappa$. But the
latter map is an isomorphism, since the same holds for
$\phi_\lambda$ and $\phi_\mu$. By Nakayama's lemma, we
deduce that $N_{\lambda\mu}=0$, therefore
$M=M_\lambda/IM_\lambda$ for any $\lambda\in\Lambda$, whence (c).

Lastly, (c)$\Rightarrow$(a) is obvious.
\end{proof}

\begin{proposition}\label{prop_loc-measur-Jacob}
Let $A$ be a locally measurable $K^+$-algebra. Suppose that
\begin{enumerate}
\alphaenu
\item
$A/\fm_KA$ is a noetherian ring.
\item
The Jacobson radical of $A$ contains $\fm_KA$.
\end{enumerate}
Then we have :
\begin{enumerate}
\item
$A\otimes_{K^+}K$ is a noetherian ring.
\item
Every finitely generated $K^+$-flat $A$-module is finitely presented.
\item
If the valuation of $K$ is discrete, $A$ is noetherian.
\end{enumerate}
\end{proposition}
\begin{proof} For any $K^+$-algebra $B$, and any ideal
$I\subset B_K:=B\otimes_{K^+}K$, let us set $I^\sat:=\Ker(B\to B_K/I)$,
and notice that $I^\sat\otimes_{K^=}K=I$. To begin with, we remark :
\begin{claim}\label{cl_saturate-measure}
Let $B$ be a measurable $K^+$-algebra, $I\subset B_K$ an ideal.
Then $I^\sat$ is a finitely generated ideal of $B$.
\end{claim}
\begin{pfclaim} Clearly $B/I^\sat$ is a $K^+$-flat finitely
generated $B$-module, hence it is finitely presented, by
proposition \ref{prop_Gruson-further}; now the claim follows
from \cite[Lemma 2.3.18(ii)]{Ga-Ra}.
\end{pfclaim}

(i): Suppose $(I_k~|~k\in\N)$ is an increasing sequence of
ideals of $A\otimes_{K^+}K$; assumption (a) implies that
there exists $n\in\N$ such that the images of $I^\sat_n$
and $I^\sat_m$ agree in $A/\fm_KA$, for every $m\geq n$.
This means that
\set\begin{equation}\label{eq_allalla}
(I^\sat_m/I^\sat_n)\otimes_{K^+}\kappa=0
\qquad
\text{for every $m\geq n$}.
\end{equation}
Next, for any geometric point $\xi$ of $\Spec\,A$, let
$A^\sh_\xi$ be the strict henselization of $A$ at $\xi$;
since the natural map $A\to A^\sh_\xi$ is flat,
$I_k\otimes_AA^\sh_\xi$ is an ideal of $A^\sh_\xi\otimes_{K^+}K$,
and clearly
$I_k^\sat\otimes_AA^\sh_\xi=(I_k\otimes_AA^\sh_\xi)^\sat$,
for every $k\in\N$. Especially $I_k^\sat\otimes_AA^\sh_\xi$
is a finitely generated ideal of $A^\sh_\xi$, by claim
\ref{cl_saturate-measure}. In view of \eqref{eq_allalla},
and assumption (b), Nakayama's lemma then says that
$(I^\sat_m/I^\sat_n)\otimes_AA^\sh_\xi=0$ for every
$m\geq n$. Since $\xi$ is arbitrary, we conclude that
$I^\sat_m/I^\sat_n=0$ for every $m\geq n$, therefore
the sequence $(I_k~|~k\in\N)$ is stationary.

(ii): Let $M$ be a $K^+$-flat and finitely generated $A$-module,
pick an $A$-linear surjection $\phi:A^{\oplus n}\to M$, and
set $N:=\Ker\,\phi$. Since $M$ is $K^+$-flat, $N/\fm_KN$ is
the kernel of $\phi\otimes_{K^+}\kappa$, and assumption
(a) implies that $N/\fm_KN$ is a finitely generated
$A/\fm_KA$-module. Hence we may find a finitely generated
$A$-submodule $N'\subset N$ such that $N=N'+\fm_KN$.
For any geometric point $\xi$ of $\Spec\,A$, define
$A^\sh_\xi$ as in the foregoing; by proposition
\ref{prop_Gruson-further} (and by
\cite[Lemma 2.3.18(ii)]{Ga-Ra}), $N\otimes_AA^\sh_\xi$
is a finitely generated $A^\sh_\xi$-module; then (b)
and Nakayama's lemma imply that
$N'\otimes_AA^\sh_\xi=N\otimes_AA^\sh_\xi$. Since $\xi$
is arbitrary, it follows that $N=N'$; especially, $M$
is finitely presented, as stated.

(iii): It suffices to show that every prime ideal of
$A$ is finitely generated (\cite[Th.3.4]{Mat}). However,
let $\fp\subset A$ be such a prime ideal, and fix a generator
$t$ of $\fm_K$. Suppose first that $t\in\fp$; in that case (a)
implies that $\fp/tA$ is a finitely generated ideal of $A/tA$,
so then clearly $\fp$ is finitely generated as well.
Next, in case $t\notin\fp$, the quotient $A/\fp$ is
a $K^+$-flat finitely generated $A$-module, hence it is
finitely presented (proposition \ref{prop_Gruson-further}),
so again $\fp$ is finitely generated
(\cite[Lemma 2.3.18(ii)]{Ga-Ra}).
\end{proof}

\begin{theorem}\label{th_finite-rank-meas}
Let $A$ be a locally measurable $K^+$-algebra, $M$ a finitely
presented $A$-module, and suppose that :
\begin{enumerate}
\alphaenu
\item
The valuation of $K$ has finite rank.
\item
For every $t\in\Spec\,K^+$, the ring $A\otimes_{K^+}\kappa(t)$
is noetherian.
\end{enumerate}
Then we have :
\begin{enumerate}
\item
$M$ admits a $K^+$-flattening sequence.
\item
$\Ass\,M$ is a finite set.
\end{enumerate}
\end{theorem}
\begin{proof}(i): Set $X:=\Spec\,A$, $S:=\Spec\,K^+$, let
$f:X\to S$ be the structure morphism, and denote by $\cM$
the quasi-coherent $\cO_{\!X}$-module associated with $M$.
Also, for every quasi-coherent $\cO_{\!X}$-module $\cF$
and every $c\in K^+$, let $c_\cF:\cF\to c\cF$ be the unique
$\cO_{\!X}$-linear morphism such that $c_\cF(U)=c_{\cF(U)}$
for every affine open subset $U\subset X$ (notation of
definition \ref{def_breaks}(ii)). For every $x\in X$, set
$K^+_{f(x)}:=\cO_{\!S,f(x)}$; then $K^+_{f(x)}$ is a valuation
ring whose valuation we denote $|\cdot|_{f(x)}$, and
$\cO_{\!X,x}$ is a locally $K^+_x$-measurable algebra.
By proposition \ref{prop_flattening-strata}, the
$\cO_{\!X,x}$-module $\cM_x$ admits a $K^+_{f(x)}$-flattening
sequence $(b_{x,0},\dots,b_{x,n})$; we wish to show that
there exists an open neighborhood $U(x)$ of $x$ in $X$,
such that $(b_{x,i}\cM/b_{x,i+1}\cM)_y$ is a
$K^+/b^{-1}_{x,i}b_{x,i+1}K^+$-flat module, for every
$y\in U(x)$ and every $i=0,\dots,n-1$. To this aim, we remark :

\begin{claim}\label{cl_cesello}
Let $\cN$ be a quasi-coherent $\cO_{\!X}$-module of finite type,
$x\in X$ a point, and $b\in K^+$ such that $\cN_x$ is a flat
$K^+/bK^+$-module. Then there exists an open neighborhood
$U\subset X$ of $x$ such that the map
$$
c_{\cN,y}\otimes_{K^+}\kappa(f(y)):
\cN_y\otimes_{K^+}\kappa(f(y))\to c\cN_y\otimes_{K^+}\kappa(f(y))
$$
is an isomorphism for every $y\in U$ and every $c\in K^+$ with
$\log|c|_{f(y)}<\log|b|_{f(y)}$.
\end{claim}
\begin{pfclaim} Let $z\in X(x)$ be any point, and set $t:=f(z)$;
under our assumptions, $\cN_z$ is a flat $K^+/bK^+$-module.
Therefore, the sequence of $\cO_{\!X,z}$-modules
$$
0\to c^{-1}b\cN_z\xrightarrow{\ j\ }\cN_z\xrightarrow{\ c_{\cN,z}\ }
c\cN_z\to 0
$$
is exact, for every $c\in K^+$ with $\log|c|<\log|b|$,
and {\em a fortiori\/}, whenever $\log|c|_t<\log|b|_t$.
If the latter inequality holds, the map $j\otimes_{K^+}\kappa(t)$
vanishes, hence $c_{\cN,z}\otimes_{K^+}\kappa(t)$ is injective.

Let now $y\in X$ be any point; set $u:=f(y)$, pick $c\in K^+$
with $\log|c|_u<\log|b|_u$, set
$\cK(u):=\Ker(c_\cN\otimes_{K^+}\kappa(u))$, and suppose that
$\cK(u)_y\neq 0$. The foregoing shows that in this case
$y\notin X(x)$. On the other hand, by assumption $f^{-1}(u)$
is a noetherian scheme, and $\cN(u):=\cN\otimes_{K^+}\kappa(u)$
is a coherent $\cO_{\!f^{-1}(u)}$-module, therefore
$\Ass\,\cN(u)$ is a finite set (\cite[Th.6.5(i)]{Mat}).
In view of proposition \ref{prop_misc-Ass}(ii), we conclude that
$\Ass\,\cK(u)$ is contained in the finite set
$\Ass\,\cN(u)\!\setminus\!X(x)$. Let $Z$ be the topological closure
in $X$ of the (finite) set
$\bigcup_{u\in S}\Ass\,\cN(u)\!\setminus\!X(x)$; taking into account
lemma \ref{lem_Ass-Supp}(ii,iii), we see that $U:=X\!\setminus\!Z$
will do.
\end{pfclaim}

Fix $x\in X$ and $i\leq n$; taking $\cN:=b_{x,i}\cM/b_{x,i+1}\cM$
and $b:=b^{-1}_{x,i}b_{x,i+1}$ in claim \ref{cl_cesello}, and
invoking the criterion of corollary \ref{cor_flatness-mod-I},
we obtain an open neighborhood $U_i$ of $x$ in $X$ such that
$\cN_y$ is $K^+/bK^+$-flat for every $y\in U_i$. Clearly the
subset $U(x):=U_1\cap\cdots\cap U_n$ fulfills the sought condition.
Next, pick finitely many points $x_1,\dots,x_k\in X$ and
corresponding $K^+$-flattening sequences $\underline b{}_i$ for
$\cM_{x_i}$, for each $i=1,\dots,k$, such that
$U(x_1)\cup\cdots\cup U(x_k)=X$; after reordering, the sequence
$(\underline b{}_1,\dots,\underline b{}_k)$ becomes
$K^+$-flattening for $M$ (see remark \ref{rem_breaks}(iii)).

(ii): Let $(b_0,\dots,b_n)$ be a $K^+$-flattening sequence for
$M$; in view of proposition \ref{prop_misc-Ass}(ii), it suffices
to prove that $\Ass\,b_iM/b_{i+1}M$ is a finite set, for every
$i=0,\dots,n-1$. Taking into account remark
\ref{rem_local-measura}(iv), we are then reduced to showing
\begin{claim}
Let $\cN$ be a quasi-coherent $\cO_{\!X}$-module, and $b\in K^+$
any element such that $\cN_x$ is a $K^+/bK^+$-flat and finitely
presented $\cO_{\!X,x}$-module for every $x\in X$. Then $\Ass\,\cN$
is finite.
\end{claim}
\begin{pfclaim}[] For given $x\in X$ and any geometric point
$\xi$ of $X$ localized at $x$, let $A^\sh$ be the strict
henselization of $A$ at $\xi$, and $B$ a local and essentially
finitely presented $K^+$-algebra with a local and ind-\'etale map
$B\to A^\sh$. We may assume that $\cN_x\otimes_AA^\sh$ descends
to a finitely presented $K^+/bK^+$-flat $B$-module $N_B$
(\cite[Ch.IV, Cor.11.2.6.1(ii)]{EGAIV-3}). Denote by $\bar x$
(resp. by $x_B$) the closed point of $X(\xi)=\Spec\,A^\sh$ (resp.
of $\Spec\,B$); in light of corollary \ref{cor_f-flat-invariance}
we have
\set\begin{equation}\label{eq_apply-rep}
x\in\Ass\,\cN\Leftrightarrow x\in\Ass_{\cO_{\!X,x}}\cN_x
\Leftrightarrow\bar x\in\Ass_{A^\sh}\cN_x\otimes_AA^\sh
\Leftrightarrow x_B\in\Ass_BN_B.
\end{equation}
On the other hand, it is easily seen that $\Ass_{K^+}K^+/bK^+$
consists of only one point, namely the maximal point $t$
of $\Spec\,K^+/bK^+$. From corollary \ref{cor_asso-tensor},
we deduce :
$$
x_B\in\Ass_BN_B\Leftrightarrow
x_B\in\Ass_BN_B\otimes_{K^+}\kappa(t)
$$
and by applying again repeatedly corollary \ref{cor_f-flat-invariance}
as in \eqref{eq_apply-rep}, we see that
$$
x_B\in\Ass_BN_B\otimes_{K^+}\kappa(t)\Leftrightarrow
x\in\Ass\,\cN\otimes_{K^+}\kappa(t).
$$
Summing up, we conclude that $\Ass\cN=\Ass\cN\otimes_{K^+}\kappa(t)$,
and the latter is a finite set, by \cite[Th.6.5(i)]{Mat}.
\end{pfclaim}
\end{proof}

\begin{corollary}\label{cor_absolutely-flat}
Let $(K,|\cdot|)$ be a valued field, and $A$ a locally measurable
$K^+$-algebra fulfilling conditions {\em (a)} and {\em (b)} of theorem
{\em\ref{th_finite-rank-meas}}. Let also $M$ be a finitely generated
$A$-module, and $\fp\subset A$ any prime ideal such that $M_\fp$
is a finitely presented $A_\fp$-module. We have:
\begin{enumerate}
\item
There exists $f\in A\!\setminus\!\fp$ such that $M_f$ is a
finitely presented $A_f$-module.
\item
$A$ is coherent.
\end{enumerate}
\end{corollary}
\begin{proof}(i): We may find a finitely presented $A$-module
$M'$ with an $A$-linear surjection $\phi:M'\to M$ such that
$\phi_\fp$ is an isomorphism. Set $M'':=\Ker\,\phi$; it
follows that $\Ass\,M''\subset\Ass\,M'\setminus\Spec\,A_\fp$.
However, $\Ass\,M'$ is a finite set (theorem
\ref{th_finite-rank-meas}(ii)); therefore the support of
$\Ker\,\phi$ is contained in a closed subset of $X:=\Spec\,A$
that does not contain $\fp$, and the assertion follows.

(ii): Let $M$ be a finitely presented $A$-module, $M'\subset M$
a finitely generated submodule, and let $\cM'$ be the
quasi-coherent $\cO_{\!X}$-module associated with $M'$. By remark
\ref{rem_local-measura}(iv), we know that $\cM'_x$ is a
finitely presented $\cO_{\!X,x}$-module for every $x\in X$;
from (i), it follows that $\cM'$ is a finitely presented
$\cO_{\!X}$-module, and the assertion follows easily.
\end{proof}

In view of theorem \ref{th_finite-rank-meas} and corollary
\ref{cor_absolutely-flat}, it is interesting to have criteria
ensuring that the fibres over $\Spec\,K$ of the spectrum of a
locally measurable $K^+$-algebra are noetherian. We present two
results in this direction. To state them, let us make first the
following :

\begin{definition} Let $f:X\to Y$ be a morphism of schemes.
\begin{enumerate}
\item
We say that $f$ is {\em absolutely flat\/} if the following
holds. For every geometric point $\xi$ of $X$, the induced
morphism $f_\xi:X(\xi)\to Y(f(\xi))$ is an isomorphism.
\item
We say that $f$ {\em has finite fibres\/} if the set
$f^{-1}(y)$ is finite, for every $y\in Y$.
\item
We say that $Y$ {\em admits a geometrically unibranch
stratification\/} if every irreducible closed subset $Z$
of\/ $Y$ contains a subset $U\neq\emptyset$, open in $Z$,
and geometrically unibranch.
\item
We say that a ring homomorphism $\phi:A\to B$ is
{\em absolutely flat\/} if $\Spec\,\phi$ is absolutely flat.
\end{enumerate}
\end{definition}

\begin{remark}\label{rem_geom-unibranchy}
In practice, it is often the case that a noetherian scheme
admits a geometrically unibranch stratification; for instance,
this holds (essentially by definition) for the spectrum of a
quasi-excellent noetherian ring (see definition \ref{def_excelll}).
\end{remark}

\begin{lemma}\label{lem_old-cotang}
Let $f:X\to Y$ be an absolutely flat morphism of schemes.
Then $\L_{X/Y}\simeq 0$ in $\sD(\cO_{\!X}\Mod)$.
\end{lemma}
\begin{proof} We easily reduce to the case where both $X$
and $Y$ are local schemes. Let $\xi$ be a geometric point
of $X$ localized at the closed point; then
$\L_{X(\xi)/X}\simeq 0$ in $\sD(\cO_{\!X(\xi)}\Mod)$
(\cite[Th.2.5.37]{Ga-Ra}). Hence, by transitivity
(\cite[Th.2.5.33]{Ga-Ra}), we are reduced to showing that
$\L_{X(\xi)/Y}\simeq 0$ in $\sD(\cO_{\!X(\xi)}\Mod)$.
The latter assertion assertion holds, again by virtue of
\cite[Th.2.5.37]{Ga-Ra}, since $f_\xi$ is an isomorphism.
\end{proof}

\begin{proposition}\label{prop_hens-or-norm}
Let $\phi:A\to B$ be a local and absolutely flat morphism
of local rings. Assume that either one of the following
two conditions holds :
\begin{enumerate}
\alphaenu
\item
$A$ is a normal local domain.
\item
$B$ is a henselian local ring.
\end{enumerate}
Then $\phi$ is an ind-\'etale ring homomorphism.
\end{proposition}
\begin{proof} Pick a geometric point $\xi$ of $\Spec\,B$
localized at the closed point, and let $\xi'$ denote the
image of $\xi$ in $\Spec\,A$. Let also $A^\sh$ (resp.
$B^\sh$) be the strict henselization of $A$ at $\xi'$ (resp.
of $B$ at $\xi$). Notice first that $A$ is a normal local
domain if and only if the same holds for $B$
(\cite[Ch.IV, Prop.18.8.12(i)]{EGA4}). In case (a) holds,
let $F_A$ (resp. $F_B$) be the field of fractions of $A$
(resp. $B$). We remark :

\begin{claim}\label{cl_tricky-bast}
Let $R$ be any normal local domain, and $R^\sh$ the strict
henselization of $R$ at a geometric point localized at the
closed point of $\Spec\,R$. We have :
\begin{enumerate}
\item
$R=R^\sh\cap\Frac\,R$ (where the intersection takes place
in $\Frac\,R^\sh$).
\item
For every field extension $F$ of $\Frac\,R$ contained
in $\Frac\,R^\sh$, the $R$-algebra $F\cap R^\sh$ is ind-\'etale.
\end{enumerate}
\end{claim}
\begin{pfclaim}(i): More generally, let $C\to D$ be any
faithfully flat ring homomorphism, and denote by $\Frac\,C$
(resp. $\Frac\,D$) the total ring of fractions of $C$ (resp.
of $D$); then it is easily seen that $C=D\cap\Frac\,C$, where
the intersection takes place in $\Frac\,D$ : the proof shall
be left as an exercise for the reader.

(ii): We easily reduce to the case where $F$ is a finite extension
of $F_R:=\Frac\,R$. Then, let $E$ be a finite Galois extension of
$F_R$, contained in a fixed separable closure $F^\sep$ of
$F_R^\sh:=\Frac\,R^\sh$, and containing $F$; also, denote by $R^\nu$
(resp. $R_E^\nu$) the integral closure of $R$ in $F^\sep$ (resp. in
$E$). Recall that there exists a unique maximal ideal
$\fp^\sep\subset R^\nu$ lying over the maximal ideal of $R$, such
that $R^\sh=(R^\nu)_{\fp^\sep}^I$, where $I\subset\Gal(F^\sep/F_R)$
is the inertia subgroup associated with $\fp^\sep$
(\cite[Ch.X, \S2, Th.2]{Ray}); therefore, $F_R^\sh=(F^\sep)^I$.
On the other hand, $F=E^H$ for a subgroup $H\subset\Gal(E/F_R)$,
and recall that the natural surjection
$\Gal(F^\sep/F_R)\to\Gal(E/F_R)$ maps $I$ onto the inertia subgroup
$I_E$ associated with the maximal ideal
$\fp_E:=\fp^\sep\cap E\subset R^\nu_E$. It follows that $H$
contains $I_E$. Set $\fp_F:=\fp_E\cap F$; then
$R_F:=(R^\nu_E\cap F)_{\fp_F}$ is a faithfully flat and essentially
\'etale $R$-algebra (\cite[Ch.X, \S1, Th.1(1)]{Ray}). Especially,
$R^\sh$ is also a strict henselization of $R_F$, and notice that
$\Frac\,R_F=F$, so the assertion follows from (i).
\end{pfclaim}

By assumption, $\phi$ extends to an isomorphism of local domains
$A^\sh\isom B^\sh$; in view of claim \ref{cl_tricky-bast}(i), we
see that $B=F_B\cap A^\sh$, and then the proposition follows from
claim \ref{cl_tricky-bast}(ii).

Next, suppose assumption (b) holds. Then, $\phi$ extends uniquely
to a local ring homomorphism $\phi^\he:A^\he\to B$ from the
henselization of $A$ (\cite[Ch.IV, Th.18.6.6(ii)]{EGA4});
clearly $\phi^\he$ is still absolutely flat, hence we may
replace $A$ by $A^\he$ and assume that both $A$ and $B$ are
henselian. Since $\phi$ is absolutely flat, the residue
field extension $\kappa(A)\to\kappa(B)$ is separable and
algebraic; then for every finite extension $E$ of $\kappa(A)$
contained in $\kappa(B)$ there exists a finite \'etale local
$A$-algebra $A_E$, unique up to isomorphism, whose residue
field is $E$ (\cite[Ch.IV, Prop.18.5.15]{EGA4}). The natural
map $\kappa(B)\to E\otimes_{\kappa(A)}\kappa(B)$ admits
a well defined section, given by the multiplication in
$\kappa(B)$, and the latter extends to a section $s_E$ of the
induced finite \'etale ring homomorphism $B\to A_E\otimes_AB$
(\cite[Ch.IV, Th.18.5.11]{EGA4}). The composition of $s_E$ with
the natural map $A_E\to A_E\otimes_AB$ is a ring homomorphism
$A_E\to B$ that extends $\phi$. The construction is clearly
compatible with inclusion of subextensions $E\subset E'$, hence
let $A'$ be the colimit of the system $(A_E~|~E\subset\kappa(B))$;
summing up, $\phi$ extends to a local and absolutely flat
map $A'\to B$. We may therefore replace $A$ by $A'$ and
assume from start that $\kappa(A)=\kappa(B)$. In this case,
set $B':=A^\sh\otimes_AB$, and notice that
$\kappa(A^\sh)\otimes_{\kappa(A)}\kappa(B)=\kappa(A^\sh)$,
so that $B'$ is a local ring, ind-\'etale, faithfully flat
over $B$, and with separably closed residue field. Since $B$
is already henselian, it follows that $B'$ is a strict
henselization of $B$, and therefore the induced map
$A^\sh\to B'$ is an isomorphism; by faithfully flat
descent, we deduce that $\phi$ is already an isomorphism,
and the proof of the proposition is concluded.
\end{proof}

\begin{lemma}\label{lem_fp-cloimmer}
Let $A$ be a locally measurable $K^+$-algebra, and set
$X:=\Spec\,K$. Let also $f:X\to Y$ be a morphism of
$K^+$-schemes with $\Omega_{X/Y}=0$, and $Y$ a finitely
presented $K^+$-scheme. Then $f_\xi$ is a finitely presented
closed immersion, for every geometric point $\xi$ of $X$.
\end{lemma}
\begin{proof} Fix such a geometric point $\xi$; by definition,
we may find a local $K^+$-scheme $Z$ essentially finitely
presented, and an ind-\'etale morphism $g:X(\xi)\to Z$
of $K^+$-schemes.
We may also assume that the morphism $X(\xi)\to Y$ deduced
from $f$ factors through $g$, in which case the resulting
morphism $h:Z\to Y$ is essentially finitely presented.
Under the current assumptions, $\Omega_{X(\xi)/Y}=0$, and
then lemma \ref{lem_old-cotang} easily implies that
$\Omega_{Z/Y}$ vanishes as well, so $h$ is essentially
unramified, and therefore it factors as the composition
of a finitely presented closed immersion $Z\to Z'$
followed by an essentially \'etale morphism $Z'\to Y$
(\cite[Ch.IV, Cor.18.4.7]{EGA4}). The lemma is an immediate
consequence.
\end{proof}

\begin{proposition} Let $A$ be a local, henselian, and locally
measurable $K^+$-algebra, whose structure map $K^+\to A$ is
local. Then $A$ is a measurable $K^+$-algebra.
\end{proposition}
\begin{proof} Pick a local and essentally finitely presented
$K^+$-algebra $B$ and a local and ind-\'etale map $\phi:B\to A$
of $K^+$-algebras. Let $\xi$ be a geometric point of $\Spec\,A$
localized at the closed point, and denote by $\xi'$ the image
of $\xi$ in $\Spec\,B$; also, let $A^\sh$ (resp. $B^\sh$) be
the strict henselization of $A$ at $\xi$ (resp. of $B$ at $\xi'$).
Then $\Omega_{A/B}=0$, hence the induced map $\phi^\sh:B^\sh\to A^\sh$
is a finitely presented surjection (lemma \ref{lem_fp-cloimmer}),
especially $I:=\Ker\,\phi^\sh$ is a finitely generated ideal, and
it follows as well that the residue field extension
$\kappa(B)\to\kappa(A)$ is separable and algebraic. Let $B^\he$
be the henselization of $B$, and for every field extension
$E$ of $\kappa(B)$ contained in $\kappa(A)^\sep=\kappa(B)^\sep$,
let $B^\he_E$ denote the local ind-\'etale $B^\he$-algebra --
determined up to isomorphism -- whose residue field is $E$
(\cite[Ch.IV, Prop.18.5.15]{EGA4}). We may find a finite
Galois extension $E$ of $\kappa(B)$ such that $I$ descends
to a finitely generated ideal $I_E\subset B^\he_E$. For any
automorphism $\sigma\in G:=\Gal(\kappa(B)^\sep/\kappa(B))$, let
$\bar\sigma\in G_E:=\Gal(E/\kappa(B))$ be the image of $\sigma$.
Recall that $G$ (resp. $G_E$) is the group of automorphisms
of the $B^\he$-algebra $B^\sh$ (resp. $B^\he_E$). With this
notation, we have the identity
$$
\bar\sigma(I_E)B^\sh=\sigma(I)=I
\qquad
\text{for every $\sigma\in G$}
$$
and since $B^\sh$ is a faithfully flat $B^\he$-algebra, it
follows that $I_E$ is invariant under the action of $G_E$;
by Galois descent, we conclude that $I$ descends to a finitely
generated ideal $I_0\subset B^\he$. Set $C:=B^\he/I_0$,
and let $C^\sh$ denote the strict henselization of $C$
at (the unique lifting of) the geometric point $\xi'$.
Since $A$ is henselian, $\phi$ extends uniquely to a local
homomorphism $\phi^\he:B^\he\to A$
(\cite[Ch.IV, Th.18.6.6(ii)]{EGA4}), and it is easily seen
that $\phi^\he$ factors through $C$. By construction, the
resulting map $C\to A$ is absolutely flat, so the assertion
follows from proposition \ref{prop_hens-or-norm}.
\end{proof}

\begin{lemma}\label{lem_parti-Z}
Let $X$ be a noetherian scheme that admits a geometrically
unibranch stratification, and $\cF$ a coherent
$\cO_{\!X_\et}$-module. Then there exists a partition
$$
X=X_1\cup\dots\cup X_k
$$
of $X$ into finitely many disjoint irreducible locally closed
subsets, such that the following holds. For every $i=1,\dots,k$,
every geometric point $\bar x$ of $X_i$, and every generization
$\bar u$ of $\bar x$ in $X_i$, every strict specialization map
$$
s_{\bar x,\bar u}:\cF_{\bar x}\to\cF_{\bar u}
$$
is injective (see \eqref{subsec_strict-special}).
\end{lemma}
\begin{proof} Arguing by noetherian induction, it suffices to
show that every reduced and irreducible closed subscheme $W$
of $X$ contains a subset $U\neq\emptyset$ that is open in $W$,
and such that every strict specialization map $s_{\bar x,\bar u}$
with $\bar x$, $\bar u$ localized in $U$, is injective. However,
$W$ contains an open subset $U\neq\emptyset$ such that
\begin{enumerate}
\alphaenu
\item
$U$ is geometrically unibranch
\item
$U$ is affine and irreducible
\item
$\cF$ is normally flat along $W$ at every point of $U$
(see \cite[Ch.IV,\S6.10.1]{EGAIV-2}).
\end{enumerate}
Indeed, (a) holds by assumption; (b) can be easily arranged by
shrinking $U$, since $X$ is noetherian. Lastly, (c) follows from
\cite[Ch.IV, Prop.6.10.2]{EGAIV-2}. We claim that such $U$
will do. Indeed, let $i:W\to W$ be the closed immersion; set
$$
\cI:=\Ker(\cO_{\!X}\to i_*\cO_W)
\qquad
\cR^\bullet:=\bigoplus_{n\in\N}\cI^n\cF/\!\cI^{n+1}\cF
$$
and recall that condition (c) means that $\cR^\bullet_u$
is a flat $\cO_{W,u}$-module, for every $u\in U$. Now, for
a given $s_{\bar x,\bar u}$ as in the foregoing, denote by
$x$ (resp. $u$) the support of $\bar x$ (resp. of $\bar u$)
and define descending filtrations by the rule :
$$
\Fil^k\cF_{\bar x}:=\cI^k_x\cdot\cF_{\bar x}
\qquad
\Fil^k\cF_{\bar u}:=\cI^k_u\cdot\cF_{\bar u}
\qquad
\text{for every $k\in\N$.}
$$
By \cite[Th.8.9]{Mat}, both these filtrations are separated, and
obviously $s_{\bar x,\bar u}$ is a map of filtered modules; thus,
it suffices to show that the induced maps
$\gr^k\cF_{\bar x}\to\gr^k\cF_{\bar u}$
of associated graded modules are injective, for every $k\in\N$.
However, notice the natural identifications :
$$
\gr^\bullet\cF_{\bar x}=
\cR^\bullet_x\otimes_{\cO_{W,x}}\cO_{W_\et,\bar x}
\qquad
\gr^\bullet\cF_{\bar u}=
\cR^\bullet_x\otimes_{\cO_{W,x}}\cO_{W_\et,\bar u}.
$$
Then, the normal flatness condition reduces to checking that
the induced strict specialization map
$\cO_{W_\et,\bar x}\to\cO_{W_\et,\bar u}$ is injective;
the latter assertion holds by the following :

\begin{claim}
Let $W$ be a reduced, irreducible scheme, $\bar w$ a geometric
point of $W$, and suppose that $W$ is unibranch at the support
$w$ of $\bar w$. Then, for every generization $\bar u$ of $\bar w$
in $W$, every strict specialization map
$\cO_{W_\et,\bar w}\to\cO_{W_\et,\bar u}$ is injective.
\end{claim}
\begin{pfclaim}[] Let $W^\nu$ be the normalization of $W$, and
$\bar w{}^\nu$ a geometric point of $W^\nu$ whose image
in $X$ is isomorphic to $\bar w$, and denote by $w^\nu$ the
support of $\bar w{}^\nu$. The assumption on $w$ means that the
induced morphism $W^\nu(w^\nu)\to W(w)$ is integral, and the
residue field extension $\kappa(w)\to\kappa(w^\nu)$ is radicial,
hence the natural morphism of $W^\nu(\bar w{}^\nu)$-schemes
$$
W(\bar w)\times_{W(w)}W^\nu(w^\nu)\to W^\nu(\bar w{}^\nu)
$$
is an isomorphism (\cite[Ch.IV, Prop.18.8.10]{EGA4}).
Since $W^\nu(\bar w{}^\nu)$ is a normal local scheme
(\cite[Prop.18.8.12(i)]{EGA4}), it follows easily that
$W(\bar w)$ is reduced and irreducible. However, any
specialization map is the composition of a localization
map, followed by a local ind-\'etale map of local rings,
whence the claim.
\end{pfclaim}
\end{proof}

\sset\subsubsection{}\label{subsec_absolutely-flat}
Now, consider -- quite generally -- a ring homomorphism $A\to B$,
with $A$ noetherian. Set $Y:=\Spec\,A$, $X:=\Spec\,B$, and denote
by $f:X\to Y$ the associated morphism of affine schemes.

\begin{lemma}\label{lem_absolu-flat}
In the situation of \eqref{subsec_absolutely-flat}, suppose
moreover that $f$ is absolutely flat. Then the following
conditions are equivalent :
\begin{enumerate}
\alphaenu
\item
$B$ is noetherian.
\item
$f$ has finite fibres.
\end{enumerate}
\end{lemma}
\begin{proof} Let $\xi$ be any geometric point of $X$, and
$x$ (resp. $y$) the support of $\xi$ (resp. of $f(\xi)$).

(a)$\Rightarrow$(b): $X(\xi)\times_{\cO_{Y,y}}\kappa(y)$ is the
strict henselization of $X_y:=f^{-1}(y)$ (with its reduced subscheme
structure) at $\xi$. Under our assumptions, $X_y(\xi)$ is isomorphic
to $\Spec\,\kappa(y)$, and (a) implies that $X_y(x)$ is a noetherian
scheme of dimension $\leq 0$. Since $x$ is arbitrary, (b) follows.

(b)$\Rightarrow$(a): Let $I\subset B$ be any ideal, and denote
by $\cI$ the associated quasi-coherent $\cO_{\!X}$-module. We
need to show that $\cI$ is an $\cO_{\!X}$-module of finite type.
Since $f_\xi$ is an isomorphism, we may find a commutative diagram
of affine schemes
$$
\xymatrix{ X' \ar[r]^-g \ar[d]_h & Y' \ar[d] \\
           X \ar[r]^-f & Y
}$$
whose left (resp. right) vertical arrow is an \'etale neighborhood
of $\xi$ (resp. of $f(\xi)$), and a quasi-coherent ideal
$\cJ\subset\cO_{\!Y'}$ such that $g^*\!\!\cJ\subset h^*\!\cI$
and such that $(g^*\!\!\cJ)_{x'}=(h^*\!\cI)_{x'}$ for some
$x'\in h^{-1}(x)$ (notice that $g$ is a flat morphism). Set
$\cM:=\cO_{Y'}/\!\cJ$, let $z'\in X'$ an arbitrary point, and
set $y':=g(z')$; in light of corollary \ref{cor_f-flat-invariance},
we have
$$
z'\in\Ass\,g^*\cM\Leftrightarrow z'\in\Ass\,g^*\cM_{|X'(x')}
\Leftrightarrow y'\in\Ass\,\cM_{|Y'(y')}\Leftrightarrow
y'\in\Ass\,\cM.
$$
In other words, $\Ass\,g^*\cM=g^{-1}\Ass\,\cM$. However,
$\Ass\,\cM$ is finite (\cite[Th.6.5(i)]{Mat}), hence the same
holds for $\Ass\,g^*\cM$, in view of (b). Now, notice that
$\Ass\,h^*\!\cI/g^*\!\!\cJ\subset\Ass\,g^*\cM\setminus X'(x')$
(proposition \ref{prop_misc-Ass}(ii)). It follows that there
exists an open neighborhood $U$ of $x'$ in $X'$ such that
$(h^*\!\cI/g^*\!\!\cJ)_{|U}=0$ (lemma \ref{lem_Ass-Supp}(iii)),
{\em i.e.} $(g^*\!\!\cJ)_{|U}=h^*\!\cI_{|U}$; especially,
$h^*\!\cI_{|U}$ is an $\cO_{\!U}$-module of finite type.
Since $h$ is an open map, we deduce that $h(U)$ is an open
neighborhood of $x$ in $X$, and $\cI_{|h(U)}$ is an
$\cO_{\!h(U)}$-module of finite type. Since $x$ is arbitrary,
the lemma follows.
\end{proof}

\sset\subsubsection{}\label{subsec_second-critter}
For our second criterion, keep the situation of
\eqref{subsec_absolutely-flat}, and suppose additionaly that,
for every geometric point $\xi$ of $X$, the induced morphism
$f_\xi:X(\xi)\to Y(\xi)$ is a closed immersion. Under this
weaker assumption, it is not necessarily true that conditions
(a) and (b) of lemma \ref{lem_absolu-flat} are equivalent.
For instance, we have :

\begin{example} Take $A:=k[X]$, the free polynomial algebra
over a given infinite field $k$; also, let $(a_i~|~i\in\N)$
be a sequence of distinct elements of $k$. We construct an
$A$-algebra $B$, as the colimit of the inductive system
$(B_i~|~i\in\N)$ of $A$-algebras, such that 
\begin{itemize}
\item
$B_i:=k[X]\times k^{i+1}$ (the product of $k[x]$ and $i+1$
copies of $k$, in the category of rings)
\item
the structure map $A\to B_i$ is the unique map of $k$-algebras
given by the rule : $X\mapsto(X,a_0,\dots,a_i)$
\item
the transition maps $B_i\to B_{i+1}$ are given by the rule :
$X\mapsto(X,0,\dots,0,a_{i+1})$ and $(0,e_i)\mapsto(0,e_i,0)$
for $i=0,\dots,i$. (Here $e_0,\dots,e_i$ is the standard basis
of the $k$-vector space $k^{i+1}$.)
\end{itemize}
Then one can check that $X:=\Spec\,B=\Spec\,k[X]\cup\N$, and
$\N$ is an open subset of $X$ with the discrete topology.
Moreover, the induced map $\Spec\,B\to Y:=\Spec\,A$ restricts
to the continuous map $\N\to Y$ given by the rule $i\mapsto\fp_i$,
for every $i\in\N$, where $\fp_i$ is the prime ideal generated
by $X-a_i$. It follows easily that the condition of
\eqref{subsec_second-critter} is fulfilled; nevertheless,
clearly $X$ has infinitely many maximal points, hence its
underlying topological space is not noetherian, and
{\em a fortiori}, $B$ cannot be noetherian.
\end{example}

However, we have the following positive result :

\begin{proposition}\label{prop_geom-unibranchy}
In the situation of \eqref{subsec_second-critter}, suppose
additionally that $Y$ admits a geometrically unibranch
stratification. Then the following conditions are equivalent :
\begin{enumerate}
\alphaenu
\item
$B$ is noetherian.
\item
The topological space underlying $X$ is noetherian.
\item
For every geometric point $\xi$ of $X$ there exist a
neighborhood $X'$ of $\xi$ in $X_\et$, an unramified
$Y$-scheme $Y'$, and an absolutely flat morphism
$X'\to Y'$ of\/ $Y$-schemes, with finite fibres.
\end{enumerate}
\end{proposition}
\begin{proof} Obviously (a)$\Rightarrow$(b).

(c)$\Rightarrow$(a): Indeed, under assumption (c), we may find
finitely many geometric points $\xi_1,\dots,\xi_k$ of $X$,
and for every $i=1,\dots,k$, a neighborhood $X'_i$ of $\xi_i$
in $X_\et$, and a $Y$-morphism $X'_i\to Y'_i$ with the stated
properties, such that moreover, the family $(X'_i~|~i=1,\dots,k)$
is an \'etale covering of $X$. Furthermore, we may assume that
$X'_i$ and $Y'_i$ are affine for every $i\leq k$. In this case,
lemma \ref{lem_absolu-flat} shows that every $X'_i$ is noetherian,
and then the same holds for $X$.

(b)$\Rightarrow$(c): Fix a geometric point $\xi$ of $X$, and
let $B^\sh_\xi$ (resp. $A^\sh_\xi$) denote the strict henselization
of $B$ at $\xi$ (resp. of $A$ at $f(\xi)$); by assumption, we may
find a (finitely generated) ideal $I\subset A^\sh_\xi$ such that
$f_\xi$ induces an isomorphism $A^\sh_\xi/I\isom B^\sh_\xi$.
Then we may find an affine \'etale neighborhood $Y'$ of $f(\xi)$,
say $Y':=\Spec\,A'$ for some \'etale $A$-algebra $A'$, and an
ideal $I'\subset A'$ such that $I'A^\sh_\xi=I$. Next, we may
find an affine \'etale neighborhood $X':=\Spec\,B'$ of $\xi$
such that $f_\xi$ extends to a morphism $g:X'\to Y'$, and we may
further suppose that the corresponding ring homomorphism $A'\to B'$
factors through $A'/I'$, so $g$ factors through a morphism
$$
h:X'\to Z:=\Spec\,A'/I'
$$
and the closed immersion $Z\to Y'$. By construction, $\xi$
lifts to a geometric point $\xi'$ of $X'$, and
$h_{\xi'}:X'(\xi')\to Z(h(\xi'))$ is an isomorphism
(\cite[Ch.IV, Prop.18.8.10]{EGA4}). To conclude the proof,
it then suffices to exhibit an open subset $U\subset X'$
containing the support of $\xi'$, and such that the
restriction $U\to Z$ of $h$ is absolutely flat with finite
fibres.

\begin{claim}\label{cl_lift-geo-unib}
Let $\phi:W\to W'$ be a quasi-finite, separated, dominant
and finitely presented morphism of reduced, irreducible
schemes, and suppose that $W'$ contains a non-empty
geometrically unibranch open subset. Then the same holds
for $W$.
\end{claim}
\begin{pfclaim} After replacing $W'$ by some open subset
$U'\subset W'$, and $W$ by $\phi^{-1}U'$, may assume that
$W'$ is affine and unibranch; then we may also suppose that
$\phi$ is finite (\cite[Ch.IV, Th.8.12.6]{EGAIV-3}), in
which case $W$ is affine as well, and $\phi$ is surjective.

Let $\eta\in W$ and $\eta'\in W'$ be the respective
generic points, and denote by $E\subset\kappa(\eta)$ the maximal
subfield that is separable over $\kappa(\eta')$. We may then
find a reduced and irreducible scheme $W''$, with generic
point $\eta''$, such that $\phi$ factors as the composition
of finite surjective morphisms $\phi':W\to W''$,
$\phi'':W''\to W'$, and such that $\kappa(\eta'')=E$.
By virtue of \cite[Ch.IV, Th.8.10.5]{EGAIV-3} and
\cite[Ch.IV, Prop.17.7.8(ii)]{EGA4}, we may then replace
$W'$ by a non-empty open subset, and assume that $\phi'$
is radicial, and $\phi''$ is \'etale. Then $W''$ is geometrically
unibranch (\cite[Ch.IV, Prop.6.15.10]{EGAIV-2}); hence, we may
replace $W'$ by $W''$, and reduce to the case where $\phi$
is radicial. Let $p$ be the characteristic of $\kappa(\eta')$;
if $p=0$, $\phi$ is birational, in which case the assertion
follows from \cite[Prop.6.15.5(ii)]{EGAIV-2}. In case $p>0$,
write $W=\Spec\,C$, $W'=\Spec\,C'$; the induced ring homomorphism
$C'\to C$ is finite and injective, and we have $C^{p^n}\subset C'$
for $n\in\N$ large enough. Denote by $C^\nu$ (resp. $C'{}^\nu$)
the normalization of the domain $C$ (resp. of $C'$); it
follows easily that $(C^\nu)^{p^n}\subset C'{}^\nu$, so
the morphism $\Spec\,C^\nu\to\Spec\,C'{}^\nu$ is radicial.
On the other hand, since $W'$ is geometrically unibranch, the
normalization morphism $\Spec\,C'{}^\nu\to W'$ is radicial
(\cite[Ch.0, Lemme 23.2.2]{EGAIV}) and therefore the
normalization map $\Spec\,C^\nu\to W$ is radicial as well
(\cite[Ch.IV, Lemme 6.15.3.1(i)]{EGAIV-2}). Then $W$ is
geometrically unibranch, again by \cite[Ch.0, Lemme 23.2.2]{EGAIV}.
\end{pfclaim}

From claim \ref{cl_lift-geo-unib} and our assumption on $Y$,
it follows easily that $Z$ admits a geometrically unibranch
stratification. The morphism $h$ induces as usual a morphism
of \'etale topoi
\set\begin{equation}\label{eq_topoi-chevreuil}
\xymatrix{ (X')^\sim_\et \ar@<.5ex>[r]^-{h_*} &
           \ar@<.5ex>[l]^-{h^*} Z^\sim_\et
}\end{equation}
as well as a morphism $h^\natural:h^*\cO_{\!Z_\et}\to\cO_{\!X'_\et}$
of $(X')^\sim_\et$-rings. By construction, for every geometric
point $\tau$ of $X'$, the induced map on stalks
$h^\natural_\tau:\cO_{\!Z_\et,h(\tau)}\to\cO_{\!X'_\et,\tau}$
is surjective, and $h^\natural_{\xi'}$ is a bijection.
It follows easily that $h^\natural_\tau$ is also bijective
for every generization $\tau$ of $\xi'$. Now, choose a partition
$Z=Z_1\cup\cdots\cup Z_k$ as in lemma \ref{lem_parti-Z}
(with $\cF:=\cO_{\!Z_\et}$), and for given $i\leq k$, suppose
that $\tau$ and $\eta$ are two geometric points of $h^{-1}Z_i$,
with $\eta$ a generization of $\tau$. The choice of a strict
specialization morphism $X'(\eta)\to X'(\tau)$ yields a
commutative diagram 
\set\begin{equation}\label{eq_chevre-commut}
{\diagram
\cO_{\!Z_\et,h(\tau)} \ar[r]^-{h^\natural_\tau} \ar[d] &
\cO_{\!X'_\et,\tau} \ar[d] \\
\cO_{\!Z_\et,h(\eta)} \ar[r]^-{h^\natural_\eta} &
\cO_{\!X'_\et,\eta}
\enddiagram}
\end{equation}
whose vertical arrows are strict specialization maps
(see remark \ref{rem_down-stri-spec}(i)); in light
of lemma \ref{lem_parti-Z}, we deduce that $h^\natural_\tau$
is injective, whenever the same holds for $h^\natural_\eta$.
Now, let $x'\in X'$ be the support of $\xi'$, and $\Sigma_i$
the set of maximal points of $h^{-1}Z_i$. Condition (b) and
proposition \ref{prop_ohm} imply that $h^{-1}Z_i$ is a noetherian
topological space, hence $\Sigma_i$ is a finite set; it follows
that the topological closure $W$ of
$\bigcup_{i=1}^k\Sigma_i\setminus X'(x')$ in $X'$ is a closed
subset that does not contain $x'$; by construction, $h^\natural$
restricts to a monomorphism on $U:=X'\setminus W$, {\em i.e.}
the restriction $h_U:U\to Z$ of $h$ is absolutely flat, as
required. It also follows that the fibres of $h_U$ are noetherian
topological spaces of dimension zero (cp. the proof of lemma
\ref{lem_absolu-flat}), hence they are finite, and the proof
is complete. 
\end{proof}

\begin{lemma}\label{lem_partitions}
Suppose that the valuation of $K$ has finite rank, and let
$X$ be a finitely presented $K^+$-scheme, $\cF$ a coherent
$\cO_{\!X_\et}$-module. Then there exists a partition
$$
X=X_1\cup\dots\cup X_k
$$
of $X$ into finitely many disjoint irreducible locally closed
subsets, such that the following holds. For every $i=1,\dots,k$,
every geometric point $\bar x$ of $X_i$, and every generization
$\bar u$ of $\bar x$ in $X_i$, every strict specialization map
$$
s_{\bar x,\bar u}:\cF_{\bar x}\to\cF_{\bar u}
$$
is injective (see \eqref{subsec_strict-special}).
\end{lemma}
\begin{proof} We easily reduce to the case where $X$ is affine,
say $X=\Spec\,A$, and then $\cF$ is the coherent
$\cO_{\!X_\et}$-module arising from a finitely presented
$A$-module $M$. By theorem \ref{th_finite-rank-meas}(i),
$M$ admits a $K^+$-flattening sequence $(b_0,\dots,b_n)$.
Then we are further reduced to showing the assertion for the
subquotients $b_iM/b_{i+1}M$ (that are finitely presented,
by corollary \ref{cor_absolutely-flat}(ii)). So, we may
assume from start that $f:X\to S_0:=\Spec\,K^+/bK^+$ is
a finitely presented morphism for some $b\in K^+$, and
$\cF$ is $f$-flat. For every $t\in S_0$, let
$$
i_t:X_t:=f^{-1}(t)\to X
$$
be the locally closed immersion; since $X_t$ is an excellent
noetherian scheme, we may apply lemma \ref{lem_parti-Z} and
remark \ref{rem_geom-unibranchy} to produce a partition
$X_t=X_{t,1}\cup\cdots\cup X_{t,k}$ by finitely many disjoint
irreducible locally closed subsets such that, for every
$i=1,\dots,k$, every geometric point $\bar x$ of $X_{t,i}$
and every generization $\bar u$ of $\bar x$ in $X_{t,i}$,
every specialization map
$$
(i_t^*\cF)_{\bar x}\to(i^*_t\cF)_{\bar u}
$$
is injective. Since $|S_0|$ is a finite set, the lemma
will then follow from :

\begin{claim} Let $g:Y\to T$ be any finitely presented morphism
of schemes, $\cG$ a finitely presented, quasi-coherent $g$-flat
$\cO_Y$-module, and $t\in T$. Let also $\bar y,\bar u$ be two
geometric points of $g^{-1}(t)$, such that $\bar u$ is a generization
of $\bar y$. Let $s_{\bar y,\bar u}:\cG_{\bar y}\to\cG_{\bar u}$ be a
strict specialization map, and suppose that
$s_{\bar y,\bar u}\otimes_{\cO_{T,t}}\kappa(t)$ is injective.
Then the same holds for $s_{\bar y,\bar u}$.
\end{claim}
\begin{pfclaim}[] Set $Y':=Y(\bar y)$, denote by $j:Y'\to Y$
the natural morphism, and set $\cG':=j^*\cG$. The map
$s_{\bar y,\bar u}$ is deduced from a morphism $Y(\bar u)\to Y'$
of $Y$-schemes; the latter factors through a faithfully flat
morphism $Y(\bar u)\to Y'(u)$, where $u\in Y'$ is the image
of the closed point of $Y(\bar u)$. Hence, $s_{\bar y,\bar u}$
is the composition of the specialization map
$s':\cG'_{\bar y}\to\cG'_u$, and the injective map
$\cG'_u\to\cG'_{\bar u}$. Our assumption implies that
$s'\otimes_{\cO_{T,t}}\kappa(t)$ is injective, and it suffices
to show that the same holds for $s'$. However, $Y'$ is the
limit of a cofiltered system
$(j_\lambda:Y_\lambda\to Y~|~\lambda\in\Lambda)$
of local, essentially \'etale $Y$-schemes. Write
$\cG_\lambda:=j_\lambda^*\cG$ for every $\lambda\in\Lambda$,
and notice that the transition morphisms $Y_\lambda\to Y_\mu$
are faithfully flat, for every $\lambda\geq\mu$; it follows
that $\cG'_{\bar y}$ is the filtered union of the system of
modules
$(G_\lambda:=\Gamma(Y_\lambda,\cG_\lambda)~|~\lambda\in\Lambda)$
(proposition \ref{prop_dir-im-and-colim}(i)); likewise,
$\cG'_{\bar y}\otimes_{\cO_{T,t}}\kappa(t)$ is the filtered
union of the submodules
$(G_\lambda\otimes_{\cO_{T,t}}\kappa(t)~|~\lambda\in\Lambda)$.
We are then reduced to checking that all the restrictions
$s_\lambda:G_\lambda\to\cG'_u$ of $s'$ are injective, and we
know already that $s_\lambda\otimes_{\cO_{T,t}}\kappa(t)$ is
injective for every $\lambda\in\Lambda$. For every such $\lambda$,
let $u_\lambda\in Y_\lambda$ be the image of $u$; then $s_\lambda$
factors through the injective map $\cG_{\lambda,u_\lambda}\to\cG'_u$
and the specialization map
$s'_\lambda:G_\lambda\to\cG_{\lambda,u_\lambda}$. Consequently,
it suffices to show that $s'_\lambda$ is injective, and we know
already that the same holds for
$s'_\lambda\otimes_{\cO_{T,t}}\kappa(t)$. However,
$\cG_{\lambda,u_\lambda}$ is a localization $Q_\lambda^{-1}G_\lambda$,
for a multiplicative set
$Q_\lambda\subset\Gamma(Y_\lambda,\cO_{Y_\lambda})$, and $s'_\lambda$
is the localization map. The claim therefore boils down to the
assertion that, for every $\lambda\in\Lambda$ and every
$q\in Q_\lambda$, the endomorphism $q\cdot\one_{G_\lambda}$
is injective on $G_\lambda$, and our assumption already ensures
that $(q\cdot\one_{G_\lambda})\otimes_{\cO_{T,t}}\kappa(t)$
is injective on $G_\lambda\otimes_{\cO_{T,t}}\kappa(t)$.
However, let $g_\lambda:Y_\lambda\to T$ be the morphism induced
by $g$; by construction, $\cG_\lambda$ is a $g_\lambda$-flat
$\cO_{Y_\lambda}$-module, hence the contention follows from
\cite[Prop.11.3.7]{EGAIV-3}.
\end{pfclaim}
\end{proof}

\begin{proposition} Let $(K,|\cdot|)$ be a valued field, and
$A$ a locally measurable $K^+$-algebra fulfilling conditions
{\em (a)} and {\em (b)} of theorem {\em\ref{th_finite-rank-meas}}.
The following conditions are equivalent:
\begin{enumerate}
\alphaenu
\addenu\addenu
\item
For every geometric point $\xi$ of $X:=\Spec\,A$ there exist
a neighborhood $U$ of $\xi$ in $X_\et$, and a finitely presented
$K^+$-scheme $Z$ with an absolutely flat morphism of $K^+$-schemes
$X'\to Z$.
\item
$\Omega_{A/K^+}$ is an $A$-module of finite type.
\end{enumerate}
\end{proposition}
\begin{proof} (c)$\Rightarrow$(d) follows easily from lemma
\ref{lem_old-cotang}.

(d)$\Rightarrow$(c): Let $a_1,\dots,a_k\in A$ be a finite
system of elements such that $da_1,\dots,da_k$ generate
the $A$-module $\Omega _{A/K^+}$. We define a map of $K^+$-algebras
$A_0:=K^+[T_1,\dots,T_k]\to A$ by the rule: $T_i\mapsto a_i$
for $i=1,\dots,k$. Clearly $\Omega_{A/A_0}=0$. Let $Z_0:=\Spec\,A_0$,
and denote by $f:X\to Z_0$ the induced morphism of schemes.
Let $A^\sh_\xi$ (resp. $A^\sh_{0,\xi}$) be the strict
henselization of $A$ at $\xi$ (resp. of $A_0$ at $f(\xi)$).
According to lemma \ref{lem_fp-cloimmer}, the induced map
$A^\sh_{0,\xi}\to A^\sh_\xi$ is surjective, and its kernel
is a finitely generated ideal $I\subset A^\sh_{0,\xi}$.
In this situation, we may argue as in the proof of
proposition \ref{prop_geom-unibranchy}, to produce an
affine \'etale neighborhood $X'$ of $\xi$, a finitely
presented affine unramified $Z_0$-scheme $Z$, and a morphism
of $Z_0$-schemes $h:X'\to Z$ such that
$h_\tau:X'(\tau)\to Z(h(\tau))$ is a closed immersion for
every geometric point $\tau$ of $X'$, and $h_{\xi'}$ is an
isomorphism for some lifting $\xi'$ of $\xi$.

Then we consider the associated morphism of \'etale topoi
as in \eqref{eq_topoi-chevreuil} and the morphism
$h^\natural:h^*\cO_{\!Z_\et}\to\cO_{\!X'_\et}$
of $(X')^\sim_\et$-rings. Again, the induced map on stalks
$h^\natural_\tau$ is surjective for every geometric point
$\tau$ of $X'$, and is bijective if $\tau$ is a generization
of $\xi'$. We pick a finite partition $Z=Z_1\cup\cdots\cup Z_k$
as in lemma \ref{lem_partitions} (for $\cF:=\cO_{\!Z}$).
For any $i\leq k$, let $\tau$, $\eta$ be two geometric
points of $h^{-1}Z_i$, such that $\eta$ is a generization
of $\tau$; by considering the commutative diagram
\eqref{eq_chevre-commut}, we see again that $h^\natural_\tau$
is injective whenever the same holds for $h^\natural_\eta$.

Now, condition (b) of theorem \ref{th_finite-rank-meas}
easily implies that $h^{-1}Z_i$ is a noetherian topological
space, hence its set $\Sigma_i$ of maximal points is finite.
Again we let $x'\in X'$ be the support of $\xi'$, and $W$
the topological closure of
$\bigcup_{i=1}^k\Sigma_i\setminus X'(x')$ in $X'$, and
it is easily seen that the restriction $U\to Z$ of $h$
is absolutely flat, so (c) holds.
\end{proof}

\sset\subsubsection{}\label{subsec_end-of-digres}
Henceforth we restrict to the case where the value group
$\Gamma$ of $K$ is not discrete and of rank one.
As usual, we consider the almost structure attached to
the standard setup attached to $(K,|\cdot|)$.

\begin{definition}
In the situation of \eqref{subsec_end-of-digres}, let
$A$ be a $K^{+a}$-algebra and $M$ an $A$-module.
\begin{enumerate}
\item
We say $M$ is an {\em almost noetherian\/} $A$-module,
if every $A$-submodule of $M$ is almost finitely generated.
\item
We say that $A$ is an {\em almost noetherian\/}
$K^{+a}$-algebra, if $A$ is an almost noetherian $A$-module.
\end{enumerate}
\end{definition}

\begin{remark}\label{rem_was-breaks}
In the situation of \eqref{subsec_end-of-digres}, suppose that
$A$ is an almost noetherian $K^{+a}$-algebra. Then the same
argument as in the ``classical limit'' case shows that every
almost finitely generated $A$-module $M$ is almost noetherian.
The details shall be left to the reader.
\end{remark}

\begin{theorem}\label{th_alm-noether-meas}
Let $A$ be a locally measurable $K^+$-algebra. Suppose that
both $A\otimes_{K^+}\kappa$ and $A\otimes_{K^+}K$ are noetherian
rings. Then $A^a$ is an almost noetherian $K^{+a}$-algebra.
\end{theorem}
\begin{proof} If the valuation of $K$ is discrete, the
assertion is proposition \ref{prop_loc-measur-Jacob}(iii).
Hence, we may assume that the valuation of $K$ is not discrete.
Moreover, let $\cQ$ be the set of all locally measurable
$K^+$-algebras $B$ that are quotients of $A$, and such that
$B^a$ is not almost noetherian. We have to show that
$\cQ=\emptyset$. However, for every $B\in\cQ$ the
$\kappa$-algebra $\bar B:=B\otimes_{K^+}\kappa$ is a quotient
of the noetherian $\kappa$-algebra $\bar A:=A\otimes_{K^+}\kappa$;
it follows that the set $\bar\cQ:=\{\bar B~|~B\in\cQ\}$
admits minimal elements, if it is not empty. In the latter
case, we may then replace $A$ by any quotient $B\in\cQ$ such
that $\bar B$ is minimal in $\bar\cQ$, and therefore assume
that for every locally measurable quotient $B$ of $A$, either
$\bar B=\bar A$, or else $B^a$ is almost noetherian.

Let $I\subset A$ be any ideal; we have to show
that $I^a$ is almost finitely generated. By assumption,
the image $\bar I$ of $I$ in $\bar A$ is finitely generated,
and the same holds for the ideal $I_K:=I\otimes_{K^+}K$ of
$A_K:=A\otimes_{K^+}K$. Thus, we may find a finitely generated
subideal $I_0\subset I$ whose image in $\bar A$ agrees
with $\bar I$, and such that $I_0\otimes_{K^+}K=I_K$.
After replacing $A$ by $A/I_0$ and $I$ by $I/I_0$, we are
then reduced to the case where both $\bar I$ and $I_K$
vanish. Let $J\subset A$ denote the kernel of the localization
map $A\to A_K$, set $S:=1+\fm_KA$, and let $B:=S^{-1}A\times A_K$.
Clearly $B$ is a faithfully flat $A$-algebra, and $A/J$ is a
$K^+$-flat $A$-module; therefore $(A/J)\otimes_AB$ is a $K^+$-flat
$B$-module of finite type, and since $A_K$ is noetherian,
proposition \ref{prop_loc-measur-Jacob}(ii) implies that
$(A/J)\otimes_AB$ is finitely presented. Then $A/J$ is
finitely presented as well, and therefore $J$ is a finitely
generated ideal. We conclude that there exists $c\in\fm_K$
such that $cJ=0$. Then notice that $J\cap cA=0$ : indeed,
if $a\in J\cap cA$, we have $a=cx$ for some $x\in A$, and
$ca=0$, therefore $c^2x=0$, so $x\in J$, and consequently
$a=cx=0$. Now, fix $b\in\fm_K$; since the valuation of $K$
has rank one, and since $I\subset J$, it follows that there
exists $n\in\N$ large enough, so that
$$
I\cap b^nA=0.
$$
Let $i_0:=\max\{i\in\N~|~I\subset b^iA\}$, and set
$$
N:=b^{i_0}A/(I+b^{i_0+1}A)
\qquad
N':=b^{i_0}A/b^{i_0+1}A.
$$
Let $\phi:N'\to N$ be the natural surjection, and set
$c^*_N:=c_N\circ\phi:N'\to cN$ for every $c\in K^+$
(notation of definition \ref{def_breaks}(ii)).

\begin{claim}\label{cl_breaking-c}
There exists $c\in K^+$ with $\log|c|<\log|b|$ such that
$c^*_N\otimes_{K^+}\kappa$ is not an isomorphism.
\end{claim}
\begin{pfclaim} Suppose that the claim fails; then it is
easily seen that no $\gamma\in\log\Gamma^+$ with $\gamma<\log|b|$
breaks $N$. For every geometric point $\xi$ of $\Spec\,A$,
let $A^\sh_\xi$ denote the strict henselization of $A$ at
$\xi$; since $A^\sh_\xi$ is a flat $A$-algebra, we have a
natural identification
$$
cN\otimes_AA^\sh_\xi=c(N\otimes_AA^\sh)
$$
of $A^\sh$-modules; therefore, no $\gamma<\log|b|$ breaks
$N\otimes_AA^\sh_\xi$. By corollary \ref{cor_flatness-mod-I}, we
deduce that $N\otimes_AA^\sh_\xi$ is a $K^+/bK^+$-flat and finitely
presented $A^\sh_\xi$-module, for every geometric point $\xi$.
Hence, $C_\xi:=\Ker\,\phi\otimes_AA^\sh_\xi$ is a finitely generated
$A^\sh_\xi$-module, and $C_\xi\otimes_{K^+}\kappa=0$, for every
geometric point $\xi$. By Nakayama's lemma, it follows that
$C_\xi=0$  for every such $\xi$, so finally $\Ker\,\phi=0$, which
means that $I\subset b^{i_0+1}A$, contradicting the choice of $i_0$.
\end{pfclaim}

Let $c$ be as in claim \ref{cl_breaking-c}, and set $d:=cb^{i_0}$;
notice the natural isomorphism of $\bar A$-modules
$$
cN\otimes_{K^+}\kappa\isom
\frac{dA}{I\cap dA}\otimes_{K^+}\kappa.
$$
Let $\phi:A\to dA/(I\cap dA)$ be the composition of $d_A:A\to dA$
and the projection $dA\to dA/(I\cap dA)$ (notation of definition
\ref{def_breaks}(ii)); by construction, $\phi\otimes_{K^+}\kappa$
is not an isomorphism, hence there exists $x\in I\cap dA$ such
that the composition $\phi_x:A\to dA/xA$ of $d_A$ and the projection
$dA\to dA/xA$ induces a map $\phi_x\otimes_{K^+}\kappa$ with
non-trivial kernel. In other words, $dA/xA$ is a cyclic module
over a locally measurable quotient $B$ of $A$ such that the
projection $\bar A\to\bar B$ is not an isomorphism, so $B^a$
is almost noetherian. Set $I_0:=(I\cap dA)/xA$; then $I_0$ is
a submodule of $dA/xA$, and consequently $I_0^a$ is an almost
finitely generated $B^a$-module (remark \ref{rem_was-breaks}).
Then clearly $(I\cap dA)^a$ is an almost finitely generated
ideal of $A^a$. But by construction, $cI\subset I\cap dA$,
and $b$ annihilates $(I\cap dA)/cI$. Since $b$ is arbitrary,
this easily implies that $I^a$ is almost finitely generated,
as required.
\end{proof}

\begin{corollary} In the situation of theorem
{\em\ref{th_alm-noether-meas}}, the following holds :
\begin{enumerate}
\item
Every almost finitely generated $A^a$-module is almost
finitely presented.
\item
Every flat almost finitely generated $A^a$-module
is almost projective of finite rank. 
\end{enumerate}
\end{corollary}
\begin{proof} Assertion (i) is an easy consequence of theorem
\ref{th_alm-noether-meas} and remark \ref{rem_was-breaks} :
the details shall be left to the reader.

(ii): Let $M$ be a flat and almost finitely generated
$A^a$-module; by (i) and \cite[Prop.2.4.18(ii)]{Ga-Ra},
$M$ is almost projective. It remains to show that
there exists $n\in\N$ such that $\Lambda^n_{A^a}M=0$,
or equivalently, that $(\Lambda^n_AM_!)^a=0$. However,
$M_!$ is a flat $A$-module, so $\Ass\,M_!\subset\Ass\,A$;
in view of theorem \ref{th_finite-rank-meas}(ii), we may
argue by induction on the cardinality $c$ of $\Ass\,A$, and
it suffices to check that, for every $\fp\in\Ass\,A$ there
exists $n\in\N$ such that $\Lambda^n_{A_\fp}(M_!)_\fp=0$.
If $c=0$, we have $A=0$, and there is nothing to prove.
Suppose that $c>0$ and the assertion is known for every
locally measurable $K^+$-algebra $B$ such that
$B\otimes_{K^+}K$ and $B\otimes_{K^+}\kappa$ are noetherian,
and such that $\Ass\,B$ has cardinality $<c$.
Especially, for a fixed $\fp\in\Ass\,A$, we can cover
$\Spec\,A_\fp\setminus\{\fp\}$ by finitely many affine
open subsets $\Spec\,B_1,\dots,\Spec\,B_k$, and then the
inductive assumption yields $n\in\N$ such that
$\Lambda^n_{B_i}(M_!\otimes_AB_i)=
\Lambda^n_{B_i}(M\otimes_{A^a}B^a_i)_!=0$ for every
$i=1,\dots,k$. In other words,
$N:=\Lambda^n_{A_\fp}(M_\fp)$ is a flat $A^a_\fp$-module
with $\Supp\,N_!\subset\{\fp\}$.

Let $A^\sh_\fp$ denote the strict henselization of $A$ at
some geometric point localized at $\fp$; on the one hand,
the $A_\fp$-algebra $A^\sh_\fp$ is faithfully flat, and
$M_!\otimes_{A}A^\sh_\fp=(M\otimes_{A^a}(A^\sh_\fp)^a)_!$.
On the other hand, exterior powers commute with arbitrary
base change; thus, we are reduced to showing :

\begin{claim} Let $B$ be a measurable $K^+$-algebra, and
$N$ a flat almost finitely generated $B^a$-module whose
support is contained in $\{s(B)\}$ (notation of
\eqref{subsec_one-more-cat}). Then $\Lambda^n_BN_!=0$ for
every sufficiently large $n\in\N$.
\end{claim}
\begin{pfclaim}[] By assumption, we may find an essentially
finitely presented $K^+$-algebra $B_0$ and an ind-\'etale
and faithfully flat map $B_0\to B$ of $K^+$-algebras.
Set $X_0:=\Spec\,B_0$ and $X:=\Spec\,B$; since $\{s(B_0)\}$
is a constructible subset of $X_0$, the natural map
$$
B\otimes_{B_0}\Gamma_{\{s(B_0)\}}\cO_{\!X_0}\to
\Gamma_{\{s(B)\}}\cO_{\!X}
$$
is an isomorphism (lemma \ref{lem_without-cohereur}(iii)).
On the other hand, there exists a finitely generated
$s(B_0)$-primary ideal $J\subset B_0$ such that the natural
map $\Gamma_{\{s(B_0)\}}\cO_{\!X_0}\to B_0/J$ is injective
(lemma \ref{lem_a&b} and theorem \ref{th_prim-dec-exist}(i)).
It follows that the natural map
$\Gamma_{\{s(B)\}}\cO_{\!X}\to B/JB$ is injective as well.
Now, $N_!$ can be written as the colimit of a filtered
system $(L_\lambda~|~\lambda\in\Lambda)$ of free $B$-modules
of finite rank; for each $\lambda\in\Lambda$, let $L^\sim_\lambda$
be the quasi-coherent $\cO_{\!X}$-module arising from $L_\lambda$,
and define likewise $N^\sim_!$; taking into account lemma
\ref{lem_flabby-Gamma_Z}(iii.b) we deduce that the natural map
$$
N_!=\Gamma_{\{s(B)\}}N^\sim_!=
\colim_{\lambda\in\Lambda}\Gamma_{\{s(B)\}}L_\lambda^\sim\to
\colim_{\lambda\in\Lambda}L_\lambda/JL_\lambda=N_!/JN_!
$$
is injective. In other words, $N$ is a $B^a/JB^a$-module.
We may then replace $B$ by $B/JB$, and assume from start
that $B$ has Krull dimension zero. In this situation, we
may find a nilpotent ideal $I\subset B$, a valuation
ring $V$ that is a measurable $K^+$-algebra, and a finitely
presented surjection $V\to B/I$ (lemma \ref{lem_find-val-ring}).
It suffices to find $n\in\N$ such that
$(\Lambda^n_BN_!)\otimes_BB/I=0$; hence, we may replace
$N$ by $N/IN$ and $B$ by $B/I$, and assume as well that
$B=V/bV$ for some $b\in V$. In this case, the assertion
follows easily from corollary \ref{cor_flat-eldivs}.
\end{pfclaim}
\end{proof}

\section{Continuous valuations and adic spaces}
In this chapter we present Huber's theory of adic spaces.

\subsection{Formal schemes}\label{sec_formal-sch}
In this section we define a category of topologically
ringed spaces that generalize the usual formal schemes from
\cite{EGAI}. This generalization was crucial for the previous
releases of this work, that followed closely Faltings's
original method for proving the almost purity theorem.
In the current release, where we have switched to the new
approach invented by Scholze, and based on his theory of
perfectoid spaces, formal schemes are still useful, but
they play a rather different role : we will employ them
to establish some foundational properties of adic spaces
(and later, of perfectoid spaces); for such purposes, the
standard adic formal schemes of \cite{EGAI} would already
suffice, but the generalization that we worked out in
the previous releases might be interesting in its own
right, and for other applications. Besides, our former
treatment included some complements on the cohomology
of formal schemes that do not appear explicitly in
\cite{EGAIII}, which will come in handy in section
\ref{subsec_adic-spaces}.

\sset\subsubsection{}\label{subsec_top-rings}
In this section we shall deal with topological rings whose
topology is linear, so that $0\in A$ admits a fundamental
system of open neighborhoods
$I_\bullet:=(I_\lambda~|~\lambda\in\Lambda)$ consisting of
ideals of $A$ (see definition \ref{def_top-ring}(v)). Then
the separated completion $A^\wedge$ of $A$ is also a topological
ring of this type : namely, its topology is linear, defined by
the system $I^\wedge_\bullet:=(I^\wedge_\lambda~|~\lambda\in\Lambda)$,
where $I^\wedge_\lambda$ denotes the topological closure of
$I_\lambda$ in $A^\wedge$, for every $\lambda\in\Lambda$.

\sset\subsubsection{}\label{subsec_top-mods}
Keep the notation of \eqref{subsec_top-rings}. We shall
consider topological $A$-modules $(M,\cT_M)$ whose topology
is $A$-linear (see definition \ref{def_top-ring}(iv)); additionally,
we shall assume that $\cT_M$ {\em is coarser than the $I_\bullet$-adic
topology}, {\em i.e.} for every open submodule $N\subset M$,
there exists $\lambda\in\Lambda$ such that $I_\lambda M\subset N$
(see remark \ref{rem_I_bullet-adic}(ii)).
(Notice that \cite[Ch.0, \S7.7.1]{EGAI} is slightly ambiguous:
it is not clear whether all the topological modules considered
there are supposed to satisfy the foregoing additional condition.)
Then, the topology of the separated completion
$(M^\wedge,\cT^\wedge_M)$ is $A^\wedge$-linear and coarser than the
$I^\wedge_\bullet$-adic topology. Let $N$ be any other topological
$A$-module of this type; notice that the topology
$\cT^\otimes_{M,N}$ on $M\otimes_AN$ (see \eqref{subsec_tensor-topol})
is also coarser than the $I_\bullet$-adic topology, so the topology
of $M\,\hat\otimes_AN$ is coarser than the $I^\wedge_\bullet$-adic
topology. We also denote
$$
\topo.\Hom_A(M,N)
$$
the $A$-module of all continuous $A$-linear maps $M\to N$.
Notice the natural identification :
\set\begin{equation}\label{eq_top.hom}
\topo.\Hom_A(M,N^\wedge)\simeq\lim_{N'\subset N}\ \colim_{M'\subset M}
\Hom_A(M/M',N/N')
\end{equation}
where $N'$ (resp. $M'$) ranges over the family of open submodules
of $N$ (resp. of $M$).

\begin{definition}\label{def_sheaves-of-spaces}
Let $X:=(\cX,J)$ be any site, and $\cC$ any other category.
\begin{enumerate}
\item
A {\em sheaf of topological spaces on $X$}, any sheaf on $X$
with values in the category $\Top$ of topological spaces (see
definition \ref{def_sheaves-with-other-values}(ii)).
Likewise, a {\em sheaf of topological groups} (resp.
{\em of topological rings})  on $X$ is a sheaf $\cA$ on $X$
with values in the category of topological groups (resp.
topological rings).
\item
Let $\cA$ be a sheaf of topological rings on $X$.
A {\em sheaf of topological $\cA$-modules\/} -- or briefly,
a {\em topological $\cA$-module\/} -- is the datum of an
$\cA$-module $\cF$, and for every $U\in\Ob(\cX)$,
a topology $\cT_U$ on $\cF(U)$, such that
$\underline\cF(U):=(\cF(U),\cT_U)$ is a topological
$\cA(U)$-module, and the rule $U\mapsto\underline\cF(U)$
defines a sheaf of topological groups on $X$.
\item
An $\cA$-linear map $\phi:\cF\to\cG$ between topological
$\cA$-modules is said to be {\em continuous\/} if it is
a morphism of sheaves of topological groups, {\em i.e.} if
$\phi_U:\cF(U)\to\cG(U)$ is a continuous map for every
$U\in\Ob(\cX)$. We denote by :
$$
\topo.\Hom_\cA(\cF,\cG)
$$
the $A$-module of all continuous $\cA$-linear morphisms
$\cF\to\cG$.
\end{enumerate}
\end{definition}

\begin{remark}\label{rem_sheaves-on-spaces}
Keep the notation of definition \ref{def_sheaves-of-spaces}.

(i)\ \
Since the forgetful functor $\Top\to\Set$ commutes with
all limits, remark \ref{rem_sheaves-with-values-in-A}(i,ii)
implies that the presheaf (of sets) underlying any sheaf of
topological spaces on $X$ is a sheaf (of sets). Likewise, the
presheaf underlying any sheaf of topological groups or topological
rings is also a sheaf of sets.

(ii)\ \ 
For any sheaf $\cA$ of topological rings on $X$, and any
topological $\cA$-modules $\cF$ and $\cG$, consider the
presheaf on $X$ defined by the rule :
\set\begin{equation}\label{eq_sheaf-top-hom}
U\mapsto\topo.\Hom_{\cA_{|U}}(\cF_{|U},\cG_{|U}).
\end{equation}
From the discussion of remark \ref{rem_sheaves-with-values-in-A}(v)
we deduce that \eqref{eq_sheaf-top-hom} is a sheaf on $X$, which we
denote by :
$$
top.\cHom_\cA(\cF,\cG).
$$
\end{remark}

\sset\subsubsection{}\label{subsec_define-space}
In the situation of \eqref{subsec_top-rings}, set
$X_\lambda:=\Spec\,A/I_\lambda$ for every $\lambda\in\Lambda$;
we define the {\em formal spectrum\/} of $A$ as the colimit of
topological spaces :
$$
\Spf\,A:=\colim_{\lambda\in\Lambda}X_{\!\lambda}.
$$
Hence, the set underlying $\Spf\,A$ is the filtered union
of the $X_\lambda$, and a subset $U\subset\Spf\,A$
is open (resp. closed) in $\Spf\,A$ if and only if
$U\cap X_\lambda$ is open (resp. closed) in $X_\lambda$
for every $\lambda\in\Lambda$.

Let $I\subset A$ be any open ideal of $A$; then $I$
contains an ideal $I_\lambda$, and therefore $\Spec\,A/I$
is a closed subset of $\Spec\,A/I_\lambda$, hence also
a closed subset of $\Spf\,A$.

\sset\subsubsection{}\label{subsec_define-sheaf}
We endow $\Spf\,A$ with a sheaf of topological rings as
follows. For every $\lambda\in\Lambda$, the structure
sheaf $\cO_{\!X_{\!\lambda}}$ carries a natural
{\em pseudo-discrete topology\/} that makes it a sheaf
of topological rings (see \cite[Ch.0, \S3.8]{EGAI}).
Let $j_\lambda:X_\lambda\to\Spf\,A$ be the closed immersion;
we set
$$
\cO_{\Spf\,A}:=
\lim_{\lambda\in\Lambda}j_{\lambda*}\cO_{\!X_{\!\lambda}}
$$
where the limit is taken in the category of sheaves of
topological rings (\cite[Ch.0, \S3.2.6]{EGAI}). By remark
\ref{rem_sheaves-with-values-in-A}(iii), it follows that :
\set\begin{equation}\label{eq_inv-limit}
\cO_{\Spf\,A}(U)=
\lim_{\lambda\in\Lambda}\cO_{X_{\!\lambda}}(U\cap X_{\!\lambda})
\end{equation}
for every open subset $U\subset\Spf\,A$. In this equality, the
right-hand side is endowed with the topology of the (projective)
limit, and the identification with the left-hand side is a
homeomorphism. Set $X:=\Spf\,A$; directly from the definitions
we get natural maps of locally ringed spaces :
\set\begin{equation}\label{eq_univ-proper}
(X_{\!\lambda},\cO_{\!X_{\!\lambda}})\stackrel{j_\lambda}{\longrightarrow}
(X,\cO_{\!X})\stackrel{i_X}{\longrightarrow}
(\Spec\,A^\wedge,\cO_{\Spec\,A^\wedge})
\qquad\text{for every $\lambda\in\Lambda$}
\end{equation}
where $A^\wedge:=\cO_{\!X}(X)$ is the (separated) completion of
$A$, and $i_X$ is given by the universal property
\cite[Ch.I, Prop.1.6.3]{EGAI-new} of the spectrum of a ring.
Clearly the composition $i_X\circ j_\lambda$ is the map of
affine schemes induced by the surjection $A^\wedge\to A/I_\lambda$.
Therefore $i_X$ identifies the set underlying $X$ with the subset
$\bigcup_{\lambda\in\Lambda}X_{\!\lambda}\subset\Spec\,A^\wedge$.
However, the topology of $X$ is usually strictly finer than
the subspace topology on the image of $i_X$. For any
$f\in A^\wedge$ we let:
$$
\fD(f):=i_X^{-1}D(f)
$$
where as usual, $D(f)\subset\Spec\,A^\wedge$ is the open subset
consisting of all prime ideals that do not contain $f$.

\sset\subsubsection{}\label{subsec_eliminate-implicit}
Let $f:A\to B$ be a continuous map of topological rings
whose topologies are linear. For every open ideal $J\subset B$,
we have an induced map $A/f^{-1}J\to B/J$, and after taking
colimits, a natural continuous map of topologically ringed spaces :
$$
\Spf\,f:(\Spf\,B,\cO_{\Spf\,B})\to(\Spf\,A,\cO_{\Spf\,A}).
$$
\begin{lemma}\label{lem_residue}
In the situation of \eqref{subsec_eliminate-implicit}, we have :
\begin{enumerate}
\item
$\Spf\,A$ is a locally and topologically ringed space.
\item
$\phi:=\Spf\,f:\Spf\,B\to\Spf\,A$ is a morphism of locally and
topologically ringed spaces; in particular, the induced map on
stalks $\cO_{\Spf\,A,\phi(x)}\to\cO_{\Spf\,B,x}$ is a local ring
homomorphism, for every $x\in\Spf\,B$.
\item
Let $I\subset A$ be any open ideal, $x\in\Spec\,A/I\subset\Spf\,A$
any point. The induced map :
\set\begin{equation}\label{eq_surge-stalk}
\cO_{\Spf\,A,x}\to\cO_{\Spec\,A/I,x}
\end{equation}
is a surjection.
\end{enumerate}
\end{lemma}
\begin{proof} (i): We need to show that the stalk $\cO_{\!X,x}$
is a local ring, for every $x\in X:=\Spf\,A$. Let
$(I_\lambda~|~\lambda\in\Lambda)$ be a cofiltered fundamental
system of open ideals. Then $x\in X_{\!\mu}:=\Spec\,A/I_\mu$ for some
$\mu\in\Lambda$; let $\kappa(x)$ be the residue field of the
stalk $\cO_{\!X_{\!\mu},x}$; there follows a ring homomorphism
\set\begin{equation}\label{eq_residue}
\cO_{\!X,x}\to(j_{\mu*}\cO_{\!X_{\!\mu}})_x=\cO_{\!X_{\!\mu},x}\to\kappa(x).
\end{equation}
Moreover, for every $\lambda\geq\mu$, the closed
immersion $X_{\!\mu}\to X_{\!\lambda}$ induces an identification of
$\kappa(x)$ with the residue field of $\cO_{\!X_{\!\lambda},x}$,
so \eqref{eq_residue} is independent of $\mu$.

Suppose now that $g\in\cO_{\!X,x}$ is mapped to a non-zero
element in $\kappa(x)$; according to \eqref{eq_inv-limit}
we may find an open subset $U\subset X$, such that $g$ is
represented as a compatible system $(g_\lambda~|~\lambda\geq\mu)$
of sections $g_\lambda\in\cO_{\!X_{\!\lambda}}(U\cap X_{\!\lambda})$.
For every $\lambda\geq\mu$, let $V_\lambda\subset X_{\!\lambda}$
denote the open subset of all $y\in X_{\!\lambda}$ such that
$g_\lambda(y)\neq 0$ in $\kappa(y)$. Then
$V_\eta\cap X_{\!\lambda}=V_\lambda$ whenever $\eta\geq\lambda\geq\mu$.
It follows that $V:=\bigcup_{\lambda\geq\mu}V_\lambda$ is
an open subset of $X$, and clearly $g$ is invertible at
every point of $V$, hence $g$ is invertible in $\cO_{\!X,x}$,
which implies the contention.

(ii): The assertion is easily reduced to the corresponding
statement for the induced morphisms of schemes :
$\Spec\,B/J\to\Spec\,A/f^{-1}J$, for any open ideal $J\subset B$.
The details shall be left to the reader.

(iii): Indeed, using \eqref{eq_univ-proper} with $X:=\Spf\,A$
and $I_\lambda:=I$, we see that the natural surjection
$\cO_{\Spec\,A^\wedge,x}\to\cO_{\Spec\,A/I,x}$ factors
through \eqref{eq_surge-stalk}.
\end{proof}

\sset\subsubsection{}\label{subsec_def-tilde}
Let $A$ and $M$ be as in \eqref{subsec_top-mods}, and fix a
(cofiltered) fundamental system $(M_\lambda~|~\lambda\in\Lambda)$
of open submodules $M_\lambda\subset M$. For every $\lambda\in\Lambda$
we may find an open ideal $I_\lambda\subset A$ such that $M/M_\lambda$
is an $A/I_\lambda$-module. Let
$j_\lambda:X_{\!\lambda}:=\Spec\,A/I_\lambda\to X:=\Spf\,A$ be
the natural closed immersion (of ringed spaces). We define
the topological $\cO_{\!X}$-module :
$$
M^\sim:=\lim_{\lambda\in\Lambda}\,j_{\lambda*}(M/M_\lambda)^\sim
$$
where, as usual, $(M/M_\lambda)^\sim$ denotes the quasi-coherent
pseudo-discrete $\cO_{\!X_{\!\lambda}}$-module associated with
$M/M_\lambda$, and the limit is formed in the category of
sheaves of topological $\cO_{\!X}$-modules
(\cite[Ch.0, \S3.2.6]{EGAI}). Thus, for every open subset
$U\subset X$ one has the identity of topological modules
\set\begin{equation}\label{eq_inv-lim-mods}
M^\sim(U)=\lim_{\lambda\in\Lambda}\,(M/M_\lambda)^\sim(U\cap X_{\!\lambda})
\end{equation}
that generalizes \eqref{eq_inv-limit}, and follows likewise
from remark \ref{rem_sheaves-with-values-in-A}(iii).

To proceed beyond these simple generalities, we need to
add further assumptions. The following definition covers
all the situations that we shall find in the sequel.

\begin{definition}\label{def_omega-form}
Let $A$ be a topological ring whose topology is linear.
\begin{enumerate}
\item
We say that $A$ is {\em $\omega$-admissible\/} if $A$ is separated
and complete, and $0\in A$ admits a countable fundamental system
of open neighborhoods.
\item
We say that an open subset $U\subset\Spf\,A$ is {\em affine\/}
if there exist an $\omega$-admissible topological ring $B$
and an isomorphism $(\Spf\,B,\cO_{\Spf\,B})\isom(U,\cO_{\Spf\,A|U})$
of topologically ringed spaces.
\item
We say that an open subset $U\subset\Spf\,A$ is
{\em truly affine\/} if $U\cap\Spec\,A/I$ is an affine open
subset of $\Spec\,A/I$ for every open ideal $I\subset A$.
\item
An {\em affine $\omega$-formal scheme\/} is a topologically
and locally ringed space $(X,\cO_{\!X})$ that is isomorphic
to the formal spectrum of an $\omega$-admissible topological ring.
\item
An {\em $\omega$-formal scheme\/} is a topologically and locally
ringed space $(X,\cO_{\!X})$ that admits an open covering
$X=\bigcup_{i\in I}U_i$ such that, for every $i\in I$, the
restriction $(U_i,\cO_{\!X|U_i})$ is an affine $\omega$-formal
scheme.
\item
A {\em morphism of $\omega$-formal schemes\/} $f:(X,\cO_{\!X})\to(Y,\cO_Y)$
is a map of topologically and locally ringed spaces, {\em i.e.}
a morphism of locally ringed spaces such that the corresponding
map $\cO_Y\to f_*\cO_{\!X}$ induces continuous ring homomorphisms
$\cO_Y(U)\to\cO_{\!X}(f^{-1}U)$, for every open subset $U\subset Y$.
\end{enumerate}
\end{definition}

\begin{remark}\label{rem_first-rem}
(i)\ \ 
Let $X$ be any $\omega$-formal scheme, and $U\subset X$ any
open subset. Using the condition of definition
\ref{def_sheaves-with-other-values}(ii), it is easily seen that :
\begin{enumerate}
\alphaenu
\item
$\cO_{\!X}(U)$ is a complete and separated topological ring.
\item
If additionally, $U=\bigcup_{i\in I}U_i$ for a countable
family $(U_i~|~i\in I)$ of open subsets, such that
$(U_i,\cO_{\!X|U_i})$ is an affine $\omega$-formal
scheme for every $i\in I$ (briefly, each $U_i$ is an
{\em affine open subset} of $X$), then $\cO_{\!X}(U)$
is $\omega$-admissible.
\end{enumerate}

(ii)\ \
In \cite{EGAI-new} one finds a notion of {\em admissible
topological ring}, and to any such ring $A$ it is assigned
an affine formal scheme $\Spf\,A$. These notions relate to
ours as follows.
Let us say that a topological ring is $c$-admissible if it is
admissible in the sense of \cite[Ch.0, D\'ef.7.1.2]{EGAI-new}
and $0\in A$ admits a countable fundamental system of open
neighborhoods. Then, an admissible topological ring $A$ is
$\omega$-admissible if and only if it is $c$-admissible.
Likewise, if $A$ is admissible, then the formal scheme
$\Spf\,A$ defined in \cite[Ch.I, D\'ef.10.1.2]{EGAI-new} is an
$\omega$-formal scheme if and only if $A$ is $c$-admissible.
Indeed, if $A$ is $c$-admissible, obviously $\Spf\,A$
coincides with the affine $\omega$-formal scheme attached
to $A$ as in \eqref{subsec_define-space} and
\eqref{subsec_define-sheaf}. For the converse, notice that
$\Spf\,A$ (in the sense of \cite{EGAI-new}) is quasi-compact,
so if it is an $\omega$-formal scheme we may cover it by
finitely many of its affine open subsets $U_1,\dots,U_n$;
say that $U_i=\Spf\,B_i$ for $i=1,\dots,n$, where
$B_1,\dots,B_n$ are $\omega$-admissible. There follows
a natural continuous map $\rho:A\to B:=\prod_{i=1}^nB_i$,
and the topology of $A$ agrees with the topology induced
by $B$ via $\rho$ (remark \ref{rem_sheaves-with-values-in-A}(i));
however, clearly $0\in B$ admits a countable fundamental
system of open neighborhoods, so the same holds for $A$.

(iii)\ \
In the same vein, if $A$ is $\omega$-admissible, then
$\Spf\,A$ as defined in \eqref{subsec_define-space} and
\eqref{subsec_define-sheaf} is a quasi-compact formal
scheme in the sense of \cite{EGAI-new} if and only if
$A$ is $c$-admissible. Indeed, the condition is obviously
sufficient; conversely, if $\Spf\,A$ is a quasi-compact
formal scheme in the sense of \cite{EGAI-new}, we may
cover it by finitely many open subsets $U_1,\dots,U_n$
such that, for every $i=1,\dots,n$ we have $U_i=\Spf\,B_i$
for some topological ring that is admissible in the sense
of \cite{EGAI-new}, and the topology of $A$ is induced
by the natural inclusion map into $B:=\prod_{i=1}^nB_i$.
For each such $B_i$, pick an open and nilpotent ideal
$I_i\subset B_i$; then $I:=\prod_{i=1}^nI_i$ is an open
and nilpotent ideal of $B$, and $I\cap A$ is open and
nilpotent in $A$, so the latter is $c$-admissible.
\end{remark}

\begin{proposition}\label{prop_truly-basis}
Suppose that $A$ is an $\omega$-admissible topological ring. Then :
\begin{enumerate}
\item
The truly affine open subsets form a basis for the topology
of $X:=\Spf\,A$.
\item
Let $I\subset A$ be an open ideal, $U\subset X$ a truly affine
open subset. Then :
\begin{enumerate}
\item
The natural map $\rho_U:A\to A_U:=\cO_{\!X}(U)$ induces an
isomorphism :
$$
(\Spf\,A_U,\cO_{\Spf\,A_U})\isom(U,\cO_{\!X|U}).
$$
\item
The topological closure $I_U$ of $IA_U$ in $A_U$
is an open ideal, and the natural map :
$\Spec\,A_U/I_U\to\Spec\,A/I$ is an open immersion.
\item
$U$ is affine.
\item
Every open covering of\/ $U$ admits a countable subcovering.
\end{enumerate}
\end{enumerate}
\end{proposition}
\begin{proof} (i): By assumption, we may find a countable
fundamental system $(I_n~|~n\in\N)$ of open ideals of $A$,
and clearly we may assume that this system is ordered under
inclusion (so that $I_n\subset I_m$ whenever $n\geq m$).
Let $x\in\Spf\,A$ be any point, and $U\subset\Spf\,A$ an
open neighborhood of $x$. Then $x\in X_n:=\Spec\,A/I_n$ for
sufficiently large $n\in\N$. We are going to exhibit, by
induction on $m\in\N$, a sequence $(f_m~|~m\geq n)$ of elements
of $A$, such that :
\set\begin{equation}\label{eq_suite}
x\in\fD(f_p)\cap X_m=\fD(f_m)\cap X_m\subset U
\qquad\text{whenever $p\geq m\geq n$}
\end{equation}
To start out, we may find $f_n\in A/I_n$ such that
$x\in\fD(f_n)\cap X_n\subset U\cap X_n$. Next, let $m\geq n$
and suppose that $f_m$ has already been found; we may
write
$$
U\cap X_{m+1}=X_{m+1}\setminus V(J)
\qquad \text{for some ideal $J\subset A/I_{m+1}$}.
$$
Let $\bar J\subset A/I_m$
and $\bar f_m\in A/I_m$ be the images of $J$ and $f_m$; then
$V(\bar J)\subset V(\bar f_m)$, hence there exists $k\in\N$
such that $\bar f{}^k_m\in\bar J$. Pick $\bar f_{m+1}\in J$
such that the image of $\bar f_{m+1}$ in $A/I_m$ agrees with
$\bar f{}^k_m$, and let $f_{m+1}\in A$ be any lifting of
$\bar f_{m+1}$; with this choice, one verifies easily that
\eqref{eq_suite} holds for $p:=m+1$. Finally, the subset :
$$
U':=\bigcup_{m\geq n}\fD(f_m)\cap X_m
$$
is an admissible affine open neighborhood of $x$ contained
in $U$.

(ii.a): For every $n,m\in\N$ with $n\geq m$, we have
a closed immersion of affine schemes: $U\cap X_m\subset U\cap X_n$;
whence induced surjections :
$$
A_{U,n}:=\cO_{\!X_n}(U\cap X_n)\to A_{U,m}:=\cO_{\!X_m}(U\cap X_m).
$$
By \cite[Ch.0, \S3.8.1]{EGAI}, $A_{U,n}$ is a discrete topological
ring, for every $n\in\N$, hence $A_U\simeq\lim_{n\in\N}A_{U,n}$
carries the linear topology that admits the fundamental
system of open ideals $(\Ker\,(A_U\to A_{U,n})~|~n\in\N)$,
especially $A_U$ is $\omega$-admissible.
It follows that the topological space underlying $\Spf\,A_U$
is $\colim_{n\in\N}(U\cap X_n)$, which is naturally identified
with $U$, under $\Spf\,\rho_U$. Likewise, let $i:U\to X$ be the
open immersion; one verifies easily ({\em e.g.} using
\eqref{eq_inv-limit}) that the natural map :
$$
i^*\lim_{n\in\N}j_{n*}\cO_{\!X_n}\to\lim_{n\in\N}i^*j_{n*}\cO_{\!X_n}
$$
is an isomorphism of topological sheaves, which implies the assertion.

(ii.b): We may assume that $I_0=I$. Since $U\cap X_n$ is affine
for every $n\in\N$, we deduce short exact sequences :
$$
\cE_n:=(0\to I\cdot\cO_{\!X_n}(X_n\cap U)\to\cO_{\!X_n}(U\cap X_n)\to
\cO_{\!X_0}(U\cap X_0)\to 0)
$$
and $\lim_{n\in\N}\cE_n$ is the exact sequence :
$$
0\to I_U\to A_U\to\cO_{\!X_0}(U\cap X_0)\to 0.
$$
Since $U\cap X_0$ is an open subset of $X_0$, both
assertions follow easily.

(ii.c): It has already been remarked that $A_U$ is
$\omega$-admissible. More precisely, for every $n\in\N$
let $(I_nA_U)^c$ be the topological closure of $I_nA_U$ in
$A_U$; then the proof of (ii.b) shows that the family of
ideals $((I_nA_U)^c~|~n\in\N)$ is a fundamental system of
open neighborhoods of $0\in A_U$. Hence the assertion follows
from (ii.a).

(ii.d): By definition, $U$ is a countable union of quasi-compact
subsets, so the assertion is immediate.
\end{proof}

\begin{corollary}\label{cor_open-of-Spf}
Let $X$ be an $\omega$-formal scheme, $U\subset X$ any open
subset. Then $(U,\cO_{\!X|U})$ is an $\omega$-formal scheme.
\qed\end{corollary}

\begin{proposition}\label{prop_univ-prop-spf}
Let $X$ be an $\omega$-formal scheme, $A$ an $\omega$-admissible
topological ring. Then the rule :
\set\begin{equation}\label{eq_nat-biject}
(f:X\to\Spf\,A)\mapsto(f^\natural:A\to\Gamma(X,\cO_{\!X}))
\end{equation}
establishes a natural bijection between the set of morphisms
of $\omega$-formal schemes $X\to\Spf\,A$ and the set of
continuous ring homomorphisms $A\to\Gamma(X,\cO_{\!X})$.
\end{proposition}
\begin{proof} We reduce easily to the case where $X=\Spf\,B$
for some $\omega$-admissible topological ring $B$. Let $f:X\to
Y:=\Spf\,A$ be a morphism of $\omega$-formal schemes; we have to
show that $f=\Spf\,f^\natural$, where $f^\natural:A\to B$ is the map
on global sections induced by the morphism of sheaves $\cO_Y\to
f^*\cO_{\!X}$ that defines $f$. Let $U\subset X$ and $V\subset Y$ be
two truly affine open subsets, such that $f(U)\subset V$, and let
likewise :
$$
f^\natural_{U,V}:A_V:=\cO_Y(V)\to B_U:=\cO_{\!X}(U)
$$
be the map induced by $f_{|U}$. Using the universal property of
\cite[Ch.I, Prop.1.6.3]{EGAI-new}, we obtain a commutative diagram
of morphisms of locally ringed spaces:
\set\begin{equation}
\diagram
\Spec\,B_U \ar[rrrr]^-{\Spec\,f^\natural_{U,V}} \ar[ddd] & & & &
\Spec\,A_V \ar[ddd] \\
& U \ar[d] \ar[rr]^-{f_{|U}} \ar[ul]^{i_U} & &
V \ar[d] \ar[ur]_{i_V} \\
& X \ar[rr]^f \ar[ld]_{i_X} & & Y \ar[dr]^{i_Y} \\
\Spec\,B \ar[rrrr]^-{\Spec\,f^\natural}& & & &\Spec\,A
\enddiagram
\end{equation}
where $i_X$, $i_Y$, $i_U$ and $i_V$ are the morphisms of
\eqref{eq_univ-proper}. Since $i_X$ and $i_Y$ are injective on the
underlying sets, it follows already that $f$ and $\Spf\,f^\natural$
induce the same continuous map of topological spaces. Let now
$J\subset B$ be any open ideal, $I\subset f^{\natural -1}(J)$ an
open ideal of $A$, and let $J_U$ (resp. $I_V$) be the topological
closure of $JB_U$ in $B_U$ (resp. of $IA_V$ in $A_V$). Since the map
$f^\natural_{U,V}$ is continuous, we derive a commutative diagram of
schemes :
$$
\xymatrix{
\Spec\,B_U/J_U \ar[r]^-\phi \ar[d]_\alpha &
\Spec\,A_V/I_V \ar[d]^\beta \\
\Spec\,B/J \ar[r]^-\psi & \Spec\,A/I }$$ where $\alpha$ and $\beta$
are open immersions, by proposition \ref{prop_truly-basis}(ii.b),
and $\phi$ (resp. $\psi$) is induced by $f^\natural_{U,V}$ (resp. by
$f^\natural$). Let $\cI$ be a fundamental system of open
neighborhoods of $0\in A$ consisting of ideals, and for
every $I\in\cI$ denote by $(IA_V)^c$ the topological
closure of $IA_V$ in $A_V$; from proposition
\ref{prop_truly-basis}(ii.b) (and its proof) it follows
as well that $A_V\simeq\lim_{I\in\cI}A_V/(IA_V)^c$. Summing up,
this shows that $f^\natural_{U,V}$ is determined by $f^\natural$,
whence the contention.
\end{proof}

\begin{corollary}\label{cor_univ-prop-spf} Let $A$ be an
$\omega$-admissible topological ring, $U\subset\Spf\,A$
an affine open subset, $V\subset U$ a truly affine open
subset of\/ $X$. Then $V$ is a truly affine open subset
of\/ $U$.
\end{corollary}
\begin{proof} Say that $U=\Spf\,B$, for some $\omega$-admissible
topological ring $B$; by proposition \ref{prop_univ-prop-spf},
the immersion $j:U\to X$ is of the form $j=\Spf\,\phi$ for
a unique map $\phi:A\to B$. Let now $J\subset B$ be any open
ideal, and set $I:=\phi^{-1}J$; there follows a commutative
diagram of locally ringed spaces :
$$
\xymatrix{
U_0:=\Spec\,B/J \ar[r] \ar[d] & X_0:=\Spec\,A/I \ar[d] \\
U \ar[r]^-j & X.
}$$
By assumption, $V\cap X_0=\Spec\,C$ for some $A/I$-algebra
$C$; it follows that :
$$
V\cap U_0=\Spec\,B/J\otimes_{A/I}C
$$
and since $J$ is arbitrary, the claim follows.
\end{proof}

\begin{remark}
(i)\ \
Let $A$ be an $\omega$-admissible topological ring, and $B$, $C$
two topological $A$-algebras whose topologies are linear (and the
structure maps $A\to B$, $A\to C$ are continuous). Recall that
$B\,\hat\otimes_AC$ represents the product of $B$ and $C$ in the
category of $\omega$-admissible $A$-algebras (see
\eqref{subsec_tensor-topol}). By standard arguments, we deduce
from proposition \ref{prop_univ-prop-spf} that
$\Spf\,B\,\hat\otimes_AC$ represents the product of $\Spf\,B$
and $\Spf\,C$ in the category of $\omega$-formal
$\Spf\,A$-schemes; especially, the category of $\omega$-formal
schemes admits arbitrary fibre products.

(ii)\ \
We shall say that a topological $A$-module $M$ is
{\em $\omega$-admissible\/} if $M$ is complete and separated,
and admits a countable fundamental system of open neighborhoods
of\ \ $0\in M$. Notice that the completion functor $N\mapsto N^\wedge$
on topological $A$-modules is not always ``right exact'', in the
following sense. Suppose $N\to N'$ is a quotient map; then the
induced map $N^\wedge\to N^{\prime\wedge}$ is not necessarily onto.
However, this is the case if $N^\wedge$ is $\omega$-admissible
(see proposition \ref{prop_replaces-Mat-Th.8.1}(v)). For such
$A$-modules, we have moreover the following :
\end{remark}

\begin{lemma}\label{lem_about-tens}
Let $A$ be a topological ring, $M$, $M'$, $N$ three
$\omega$-admissible topological $A$-modules, and
$f:M\to M'$ a quotient map. Then :
\begin{enumerate}
\item
The homomorphism
$f\,\hat\otimes_A\one_N:M\,\hat\otimes_AN\to M'\,\hat\otimes_AN$
is a quotient map.
\item
Let us endow $\Ker\,f$ with the subspace topology induced from $M$,
and suppose additionally that, for every open ideal $I\subset A$ :
\begin{enumerate}
\item
The topological closure $(IN)^c$ of $IN$ in $N$ is an open
submodule of $N$.
\item
$N/(IN)^c$ is a flat $A/I$-module.
\end{enumerate}
Then the complex :
$$
0\to (\Ker\,f)\,\hat\otimes_AN\to M\,\hat\otimes_AN\to
M'\,\hat\otimes_AN\to 0
$$
is an admissible short exact sequence of topological $A^\wedge$-modules.
\end{enumerate}
\end{lemma}
\begin{proof} By assumption, we may find an inverse system of
discrete $A$-modules, with surjective transition maps
$(M_n~|~n\in\N)$ (resp. $(N_n~|~n\in\N)$), and an isomorphism
of topological $A$-modules :
$M\simeq\lim_{n\in\N}M_n$ (resp. $N\simeq\lim_{n\in\N}N_n$);
let us define $M'_n:=M'\amalg_M M_n$ for every $n\in\N$ (the
cofibred sum of $M'$ and $M_n$ over $M$). Since
$f$ is a quotient map, $M'_n$ is a discrete $A$-module for
every $n\in\N$, and the natural map : $M'\to\lim_{n\in\N}M'_n$
is a topological isomorphism. We deduce an inverse system of
short exact sequences of discrete $A$-modules :
$$
0\to K_n\to M_n\otimes_AN_n\to M'_n\otimes_AN_n\to 0
$$
where $K_n$ is naturally a quotient of $\Ker(M_n\to M'_n)\otimes_A N_n$,
for every $n\in\N$; especially, the transition maps
$K_j\to K_i$ are surjective whenever $j\geq i$. Then
assertion (i) follows from lemma \ref{lem_crit-admissible}.

(ii): We may find open ideals $I_n\subset A$ such that
$I_nN_n=I_nM_n=0$ for every $n\in\N$, hence
$(I_nN)^c\subset\Ker(N\to N_n)$, and from (a) we deduce
that $((I_nN)^c~|~n\in\N)$ is a fundamental system of open
neighborhoods of $0\in N$. Due
to (b), we obtain short exact sequences :
$$
0\to\Ker(M_n\to M'_n)\otimes_AN/(I_nN)^c\to M_n\otimes_AN/(I_nN)^c
\to M'_n\otimes_AN/(I_nN)^c\to 0
$$
for every $n\in\N$. Then it suffices to invoke again lemma
\ref{lem_crit-admissible}.
\end{proof}

\begin{lemma}\label{lem_left-adj}
Let $A$ be an $\omega$-admissible topological ring, $M$ a topological
$A$-module, and $U\subset X:=\Spf\,A$ any affine open subset. Then the
following holds :
\begin{enumerate}
\item
There is a natural isomorphism of topological $\cO_{\!X}(U)$-modules :
$$
M^\sim(U)\simeq M\hat\otimes_A\cO_{\!X}(U).
$$
\item
If $L$ is any other topological $A$-module, the natural map :
$$
\mathrm{top}.\Hom_A(L,M^\wedge)\to
\mathrm{top}.\Hom_{\cO_{\!X}}(L^\sim,M^\sim)
\quad\phi\mapsto\phi^\sim
$$
is an isomorphism.
\item
The functor $M\mapsto M^\sim$ on topological
$A$-modules, is left adjoint to the global sections functor
$\cF\mapsto\cF(X)$, defined on the category of complete and
separated topological $\cO_{\!X}$-modules and continuous maps.
\end{enumerate}
\end{lemma}
\begin{proof}(i): First we remark that the assertion holds
whenever $U$ is a truly affine open subset of $X$; the easy
verification shall be left to the reader. For a general $U$,
set $M_U:=M\hat\otimes_A\cO_{\!X}(U)$, and denote by $M^\sim_U$
the associated $\cO_U$-module. Let $V\subset U$
be any open subset which is truly affine in $X$; by corollary
\ref{cor_univ-prop-spf}, $V$ is truly affine in $U$ as well.
We deduce natural isomorphisms:
$$
M^\sim(V)\simeq M_U\hat\otimes_{\cO_{\!X}(U)}\cO_{\!X}(V)
\simeq M^\sim_U(V)
$$
which -- in view of proposition \ref{prop_truly-basis}(i) -- amount
to a natural isomorphism of topological $\cO_{\!U}$-modules :
$(M^\sim)_{|U}\isom M^\sim_U$. Assertion (i) follows easily.

(iii): Given a continuous map $M^\sim\to\cF$, we get
a map of global sections $M=M^\sim(X)\to\cF(X)$. Conversely,
suppose $f:M\to\cF(X)$ is a given continuous map; let
$U\subset X$ be any affine open subset, and $f_U:M\to\cF(U)$
the composition of $f$ and the restriction map $\cF(X)\to\cF(U)$.
Then $f_U$ extends first -- by linearity -- to a map
$M\otimes_A\cO_{\!X}(U)\to\cF(U)$, and second -- by continuity --
to a map $M\hat\otimes_A\cO_{\!X}(U)\to\cF(U)$; the latter, in
view of (i), is the same as a map $f^\sim_U:M^\sim(U)\to\cF(U)$.
Clearly the rule $U\mapsto f^\sim_U$ thus defined is functorial
for inclusion of open subsets $U\subset U'$, whence (iii).

(ii) is a straightforward consequence of (iii).
\end{proof}

\begin{proposition}\label{prop_quasi-coh-omega}
Let $A$ be an $\omega$-admissible topological ring, $M$ an
$\omega$-admissible topological $A$-module, $N\subset M$ a
submodule, and $U\subset X:=\Spf\,A$ an affine open subset. Then:
\begin{enumerate}
\item
If we endow $M/N$ with the quotient topology, the sequence of
$\cO_{\!X}$-modules :
$$
0\to N^\sim\to M^\sim\to(M/N)^\sim\to 0
$$
is short exact.
\item
The induced sequence
$$
0\to N^\sim(U)\to M^\sim(U)\to(M/N)^\sim(U)\to 0
$$
is short exact and admissible in the sense of
\eqref{sec_Topolog-Hom-Alg}.
\end{enumerate}
\end{proposition}
\begin{proof} (i): We recall the following :

\begin{claim}\label{cl_obv-flat}
Let $j:Z':=\Spec\,R'\to Z:=\Spec\,R$ be an open immersion
of affine schemes. Then $R'$ is a flat $R$-algebra.
\end{claim}
\begin{pfclaim} The assertion can be checked on the localizations
at the prime ideals of $R'$; however, the induced maps
$\cO_{\!Z,j(z)}\to\cO_{\!Z',z}$ are isomorphisms for every
$z\in Z'$, so the claim is clear.
\end{pfclaim}

Let $V\subset X$ be any truly affine open subset. It
follows easily from proposition \ref{prop_truly-basis}(ii.b)
and claim \ref{cl_obv-flat} that $\cO_{\!X}(V)$ fulfills
conditions (a) and (b) of lemma \ref{lem_about-tens}(ii).
In light of lemma \ref{lem_left-adj}(i), we deduce that
the sequence :
$$
0\to N^\sim(V)\to M^\sim(V)\to(M/N)^\sim(V)\to 0
$$
is admissible short exact, whence the contention.

(ii): Clearly the sequence is left exact; using lemmata
\ref{lem_about-tens}(i) and \ref{lem_left-adj}(i), we see that
it is also right exact, and moreover $M^\sim(U)\to(M/N)^\sim(U)$
is a quotient map. It remains only to show that the
topology on $N^\sim(U)$ is induced from $M^\sim(U)$,
and to this aim we may assume -- thanks to the condition of
definition \ref{def_sheaves-with-other-values}(ii) -- that $U$ is
a truly affine open subset of $X$, in which case the assertion
has already been observed in the proof of (i).
\end{proof}

\begin{definition}\label{def_quasi-coh-omega}
Let $X$ be an $\omega$-formal scheme, and
$\cF$ a topological $\cO_{\!X}$-module $\cF$.

(i)\ \
We say that $\cF$ is {\em quasi-coherent\/} if there exist
a covering $\fU:=(U_i~|~i\in I)$ of $X$ consisting of affine
open subsets, and for every $i\in I$, an $\omega$-admissible
topological $\cO_{\!X}(U_i)$-module $M_i$ with an isomorphism
$\cF_{|U_i}\isom M^\sim_i$ of topological $\cO_{U_i}$-modules
(in the sense of definition \ref{def_sheaves-of-spaces}(ii,iii)).

(ii)\ \
We denote by $\cO_{\!X}\Mod_\qcoh$ the category of quasi-coherent
$\cO_{\!X}$-modules and continuous $\cO_{\!X}$-linear morphisms.

(iii)\ \
We say that a continuous morphism $f:\cF\to\cG$ of quasi-coherent
$\cO_{\!X}$-modules is a {\em quotient map\/} if the induced
map $f(U):\cF(U)\to\cG(U)$ is a quotient map of topological
$\cO_{\!X}(U)$-modules, for every affine open subset $U\subset X$.
\end{definition}

\begin{remark}
(i)\ \
The category $\cO_{\!X}\Mod_\qcoh$ is usually not abelian (see
\eqref{sec_Topolog-Hom-Alg}); more than that, the kernel
(in the category of abelian sheaves) of a continuous map
$f:\cF\to\cG$ of quasi-coherent $\cO_{\!X}$-modules may fail
to be quasi-coherent. However, using proposition
\ref{prop_quasi-coh-omega} one may show that $\Ker\,f$ is
quasi-coherent whenever $f$ is a quotient map, and in this
case $\Ker\,f$ is also the kernel of $f$ in the category
$\cO_{\!X}\Mod_\qcoh$.

(ii)\ \
Furthermore, any (continuous) morphism $f$ in $\cO_{\!X}\Mod_\qcoh$
admits a cokernel. This can be exhibited as follows. To start
out, let us define presheaves $\cI$ and $\cL$ by declaring that
$\cI(U)\subset\cG(U)$ is the topological closure of the
$\cO_{\!X}(U)$-submodule $\Img(f(U):\cF(U)\to\cG(U))$, and
$\cL(U):=\cG(U)/\cI(U)$, which we endow with its natural
quotient topology, for every open subset $U\subset X$. Now,
suppose that $V\subset U$ is an inclusion of sufficiently
small affine open subsets of $X$ (so that $\cF_{|U}$ and
$\cG_{|U}$ are of the form $M^\sim$ for some topological
$\cO_{\!X}(U)$-module $M$); since
$\cF(U)\hat\otimes_{\cO_{\!X}(U)}\cO_{\!X}(V)=\cF(V)$, we see
that the image of $\cI(U)$ in $\cI(V)$ generates a dense
$\cO_{\!X}(V)$-submodule. On the other hand, by construction
the exact sequence $\cE:=(0\to\cI(U)\to\cG(U)\to\cL(U)\to 0)$
is admissible, hence the same holds for the sequence
$\cE\hat\otimes_{\cO_{\!X}(U)}\cO_{\!X}(V)$ (lemma
\ref{lem_left-adj}(i) and proposition \ref{prop_quasi-coh-omega}(ii)).
It follows that $\cI(V)=\cI(U)\hat\otimes_{\cO_{\!X}(U)}\cO_{\!X}(V)$
and $\cL(V)=\cL(U)\hat\otimes_{\cO_{\!X}(U)}\cO_{\!X}(V)$.
Thus $\cI$ and $\cL$ are sheaves of topological $\cO_{\!X}$-modules
on the site $C$ of all sufficiently small affine open subsets of $X$,
and their sheafifications $\cI'$ and $\cL'$ are the topological
$\cO_{\!X}$-modules obtained as in \cite[Ch.0, \S3.2.1]{EGAI},
by extension of $\cI$ and $\cL$ from the site $C$ to the
whole topology of $X$. It follow that $\cI'$ and $\cL'$ are
quasi-coherent $\cO_{\!X}$-modules; then it is easy to check
that $\cL'$ is the cokernel of $f$ in the category
$\cO_{\!X}\Mod_\qcoh$ (briefly : the
{\em topological cokernel of $f$}), and shall be denoted
$$
\topo.\Coker\,f.
$$
The sheaf $\cI'$ shall be called {\em the topological closure
of the image of $f$}, and denoted
$$
\bar\Img(f).
$$
A morphism $f$ with $\topo.\Coker\,f=0$ shall be called a
{\em topological epimorphism}.

(iii)\ \
Notice also that the natural map $\Coker\,f\to\topo.\Coker\,f$
is an epimorphism of abelian sheaves, hence a continuous morphism
of quasi-coherent $\cO_{\!X}$-modules which is an epimorphism
of $\cO_{\!X}$-modules, is also a topological epimorphism.
\end{remark}

\begin{proposition}\label{prop_intrinsic-qcoh}
Let $X$ be an affine $\omega$-formal scheme, $\cF$ a quasi-coherent
$\cO_{\!X}$-module. Then $\cF(X)$ is an $\omega$-admissible
$\cO_{\!X}(X)$-module, and the natural map
$$
\cF(X)^\sim\to\cF
$$
is an isomorphism of topological $\cO_{\!X}$-modules.
\end{proposition}
\begin{proof} By assumption we may find an affine open covering
$\fU:=(U_i~|~i\in I)$ of $X$ such that, for every $i\in I$,
$\cF_{|U_i}\simeq M^\sim_i$ for some $\cO_{\!X}(U_i)$-module $M_i$.
In view of lemma \ref{lem_left-adj}(i) we may assume -- up
to replacing $\fU$ by a refinement -- that $U_i$ is a truly affine
subset of $X$ for every $i\in I$. Furthermore, we may write
$X=\bigcup_{n\in\N}X_n$ for an increasing countable family of
quasi-compact subsets; for each $n\in\N$ we may then find a finite
subset $I(n)\subset I$ such that $X_n\subset\bigcup_{i\in I(n)}U_i$,
and therefore we may replace $I$ by $\bigcup_{n\in\N}I_n$, which
allows to assume that $I$ is countable.

Next, for every $i\in I$ we may find a countable fundamental
system of open submodules $(M_{i,n}~|~n\in\N)$ of $M_i$, and
for every $n\in\N$ an open ideal $J_{i,n}\subset A_i:=\cO_{\!X}(U_i)$
such that $N_{i,n}:=M_i/M_{i,n}$ is an $A_i/J_{i,n}$-module. Let
$j_{i,n}:U_{i,n}:=\Spec\,A_i/J_{i,n}\to U_i$ be the natural closed
immersion; we may write :
\set\begin{equation}\label{eq_restrict-and-iso}
\cF_{|U_i}\simeq\lim_{n\in\N}j_{i,n*} N_{i,n}^\sim.
\end{equation}
Let also $\iota_{i,n}:U_{i,n}\to X$ be the locally closed
immersion obtained as the composition of $j_{i,n}$ and the open
immersion $j_i:U_i\to X$; we deduce natural maps of $\cO_{\!X}$-modules :
$$
\phi_{i,n}:\cF\to j_{i*}\cF_{|U_i}\to
\cG_{i,n}:=\iota_{i,n*} N_{i,n}^\sim.
$$

\begin{claim}\label{cl_cont-pseudo} (i)\ \ There exists
$m\in\N$ such that $\cG_{i,n}$ is the extension
by zero of a quasi-coherent $\cO_{\!X_m}$-module.
\begin{enumerate}
\addenu
\item
$\phi_{i,n}$ is continuous for the pseudo-discrete
topology on $\cG_{i,n}$.
\end{enumerate}
\end{claim}
\begin{pfclaim} (i): It is easy to see that $U_{i,n}\subset X_m$
for $m\in\N$ large enough; then the assertion follows from
(\cite[Ch.I, Cor.9.2.2]{EGAI}).

(ii): We need to check that the map
$\phi_{i,n,V}:\cF(V)\to\cG_{i,n}(V)$ is continuous for every
open subset $V\subset X$. However, the condition of
definition \ref{def_sheaves-with-other-values}(ii) implies that
the assertion is local on $X$, hence we may assume that $V$ is
a truly affine open subset of $X$, in which case $\cG_{i,n}(V)$
is a discrete space. We may factor $\phi_{i,n,V}$ as a composition :
$$
\cF(V)\stackrel{\alpha}{\to}\cF(V\cap U_i)
\stackrel{\beta}{\to}N^\sim_{i,n}(U_{i,n}\cap V)
$$
where the restriction map $\alpha$ is continuous, and
$\beta$ is continuous for the pseudo-discrete topology on
$N^\sim_{i,n}(U_{i,n}\cap V)$. However, $U_i\cap V$
is a truly affine open subset of $U_i$; therefore $U_{i,n}\cap V$
is quasi-compact, and the pseudo-discrete topology on $N^\sim_{i,n}$
induces the discrete topology on $N^\sim_{i,n}(U_{i,n}\cap V)$.
The claim follows.
\end{pfclaim}

For every finite subset $S\subset I\times\N$, let
$\phi_S:\cF\to\cG_S:=\prod_{(i,n)\in S}\cG_{i,n}$ be the
product of the maps $\phi_{i,n}$; according to claim
\ref{cl_cont-pseudo}(ii), $\phi_S$ is continuous for the
pseudo-discrete topology on $\cG_S$. Hence, for every $i\in I$,
the restriction $\phi_{S|U_i}:\cF_{|U_i}\to\cG_{S|U_i}$ is of the
form $f_{i,S}^\sim$, for some continuous map
$f_{i,S}:M_i\to\cG_S(U_i)$ (lemma \ref{lem_left-adj}(ii)).
It follows easily that $(\Img\,\phi_S)_{|U_i}$ is the quasi-coherent
$\cO_{U_i}$-module $(\Img\,f_{i,S})^\sim$, especially $\Img\,\phi_S$
is quasi-coherent. Notice that $\cG_S$ is already a quasi-coherent
$\cO_{\!X_m}$-module for $m\in\N$ large enough, hence the same
holds for $\Img\,\phi_S$; we may therefore find an
$\cO_{\!X_m}(X_m)$-module $G_S$ such that $G_S^\sim\simeq\Img\,\phi_S$.
Furthermore, for every other finite subset $S'\subset I\times\N$
containing $S$, we deduce a natural $\cO_{\!X}(X)$-linear surjection
$G_{S'}\to G_S$, compatible with compositions of inclusions
$S\subset S'\subset S''$. Let $(S_n~|~n\in\N)$ be a countable
increasing family of finite subsets, whose union is $I\times\N$;
we endow $G:=\lim_{n\in\N}G_{S_n}$ with the pro-discrete topology.
It is then easy to check that $G$ is an $\omega$-admissible
$\cO_{\!X}(X)$-module; furthermore, by construction we get
a unique continuous map $\phi:\cF\to G^\sim$ whose composition
with the projection onto $G_{S_n}^\sim$ agrees with $\phi_{S_n}$
for every $n\in\N$. In light of \eqref{eq_restrict-and-iso}
we see easily that $\phi_{|U_i}$ is an isomorphism of topological
$\cO_{\!U_i}$-modules for every $i\in I$, hence the same holds
for $\phi$. Then $\phi$ necessarily induces an isomorphism
$\cF(X)\isom G$.
\end{proof}

\begin{corollary}\label{cor_inv-lim-coh}
Let $X$ be an $\omega$-formal scheme, $(\cF_n~|~n\in\N)$
an inverse system of quasi-coherent $\cO_{\!X}$-modules,
whose transition maps are topological epimorphisms. Then
$\lim_{n\in\N}\cF_n$ (with its inverse limit topology)
is a quasi-coherent $\cO_{\!X}$-module.
\end{corollary}
\begin{proof} We may assume that $X$ is affine; then, for every
$n\in\N$ we have $\cF_n=M_n^\sim$ for some complete and separated
$\cO_{\!X}(X)$-module $M_n$ (by proposition \ref{prop_intrinsic-qcoh}),
and the transition maps $\cF_{n+1}\to\cF_n$ come from
corresponding continuous linear maps $f_n:M_{n+1}\to M_n$.
We choose inductively, for every $n\in\N$, a descending fundamental
system of open submodules $(N_{n,k}~|~k\in\N)$ of $M_n$, such that
$f_n(N_{n+1,k})\subset N_{n,k}$ for every $n,k\in\N$.
Set $M_{n,k}:=M_n/N_{n,k}$ for every $n,k\in\N$.
By assumption, $\topo.\Coker\,f^\sim_n=0$; hence the
induced maps $M^\sim_{n+1,k}\to M^\sim_{n,k}$ are topological
epimorphisms; since $M^\sim_{n,k}$ is pseudo-discrete (on the
closure of its support), it follows that the latter maps are
even epimorphisms (of abelian sheaves) so the corresponding
maps $M_{n+1,k}\to M_{n,k}$ are onto for every $n,k\in\N$. Thus :
$$
\lim_{n\in\N}\cF_n\simeq\lim_{n\in\N}\lim_{k\in\N}M_{n,k}^\sim
\simeq\lim_{n\in\N}M^\sim_{n,n}\simeq(\lim_{n\in\N}M_{n,n})^\sim
$$
and it is clear that the $\cO_{\!X}(X)$-module $\lim_{n\in\N}M_{n,n}$
is $\omega$-admissible.
\end{proof}

\begin{theorem}\label{th_coh-vanish}
Let $X$ be an $\omega$-formal scheme, $\cF$ a
quasi-coherent $\cO_{\!X}$-module, $\fU:=(U_i~|~i\in I)$
an open covering of $X$, such that the intersection
$U_{t_0}\cap\cdots\cap U_{t_n}$ is affine for every $n\in\N$
and every $(t_0,\dots,t_n)\in I^{n+1}$. Then :
\begin{enumerate}
\item
There is a natural isomorphism (notation of
\eqref{subsec_altern-pseudo-Leray})
$$
H^\bullet_\alt(\fU,\cF)\isom H^\bullet(X,\cF).
$$
\item
If moreover $X$ is affine, we have $H^i(X,\cF)=0$
for every $i>0$.
\end{enumerate}
\end{theorem}
\begin{proof} (The statements refer to the cohomology of the
abelian sheaf underlying $\cF$, in other words, we forget
the topology of the modules $\cF(U)$.)

(ii): Say that $X=\Spf\,A$ for an $\omega$-admissible
topological ring $A$. Let $U\subset X$ be any truly affine
open subset, and $\fU:=(U_i~|~i\in I)$ a covering of $U$
consisting of truly affine open subsets. By proposition
\ref{prop_intrinsic-qcoh} we have $\cF_{|U}\simeq M^\sim$,
where $M:=\cF(U)$.

\begin{claim}\label{cl_inv-lim-Cech}
The augmented alternating \v{C}ech complex
$C^\bullet_\mathrm{alt}(\fU,M^\sim)$ is acyclic.
\end{claim}
\begin{pfclaim} Let $(M_n~|~n\in\N)$ be a fundamental system
of neighborhoods of $0\in M$, consisting of open submodules,
and for every $n\in\N$, choose an open ideal $I_n\subset A$
such that $I_n M\subset M_n$. Set $X_n:=\Spec\,A/I_n$ and
$\fU_n:=(U_i\cap X_n~|~i\in I)$ for every $n\in\N$.
For every $n\in\N$ we may consider the augmented alternating
\v{C}ech complex $C^\bullet_\mathrm{alt}(\fU_n,(M/M_n)^\sim)$, and in
view of \eqref{eq_inv-lim-mods} we obtain a natural isomorphism
of complexes :
$$
C^\bullet_\mathrm{alt}(\fU,M^\sim)\isom
\lim_{n\in\N}C^\bullet_\mathrm{alt}(\fU_n,(M/M_n)^\sim).
$$
We may view the double complex
$C^\bullet_\mathrm{alt}(\fU_\bullet,(M/M_\bullet)^\sim)$
also as a complex of inverse systems of modules, whose term
in degree $i\in\N$ is
$C^i_\mathrm{alt}(\fU_\bullet,(M/M_\bullet)^\sim)$.
Notice that, for every $i\in\N$, all the transition maps
of this latter inverse system are surjective. Hence
(\cite[Lemma 3.5.3]{We}) :
$$
\lim_{n\in\N}{}^{\!q}\,
C^i_\mathrm{alt}(\fU_n,(M/M_n)^\sim)=0\qquad
\text{for every $q>0$.}
$$
In other words, these inverse systems are acyclic for the
inverse limit functor. It follows ({\em e.g.} by means of
\cite[Th.10.5.9]{We}) that :
$$
R\lim_{n\in\N}\,C^\bullet_\mathrm{alt}(\fU_\bullet,(M/M_\bullet)^\sim)
\simeq C^\bullet_\mathrm{alt}(\fU,M^\sim).
$$
On the other hand, the complexes
$C^\bullet_\mathrm{alt}(\fU_n,(M/M_n)^\sim)$ are acyclic for
every $n\in\N$ (\cite[Ch.III, Th.1.3.1 and Prop.1.4.1]{EGAIII}),
hence $C^\bullet_\mathrm{alt}(\fU_\bullet,(M/M_\bullet)^\sim)$
is acyclic, when viewed as a cochain complex of inverse systems
of modules. The claim follows.
\end{pfclaim}

Assertion (ii) follows from claim \ref{cl_inv-lim-Cech}, proposition
\ref{prop_truly-basis}(i) and theorem \ref{th_Cartan}(i).

Assertion (i) follows from (ii) and corollary \ref{cor_Leray}(ii).
\end{proof}

\begin{corollary}\label{cor_inverse-lim}
In the situation of \eqref{subsec_def-tilde}, suppose
that $A$ and $M$ are $\omega$-admissible. Then:
$$
\lim_{\lambda\in\Lambda}{}^{\!q}\,j_{\lambda*}(M/M_\lambda)^\sim=0
\qquad\text{for every $q>0$}.
$$
\end{corollary}
\begin{proof} For every truly affine open subset $U\subset X$,
we have a topos $U^\Lambda$ defined as in \cite[\S7.3.4]{Ga-Ra},
and the cofiltered system
$\cM:=(j_{\lambda*}(M/M_\lambda)^\sim~|~\lambda\in\Lambda)$
defines an abelian sheaf on $U^\Lambda$. According to {\em loc.cit.}
there are two spectral sequences :
$$
\begin{aligned}
E^{pq}_2 :=\, & R^p\Gamma(U,
\lim_{\lambda\in\Lambda}{}^{\!q}\,j_{\lambda*}(M/M_\lambda)^\sim)
\Rightarrow H^{p+q}(U^\Lambda,\cM) \\
F^{pq}_2 :=\, &
\lim_{\lambda\in\Lambda}{}^{\!q}\,R^p\Gamma(U,j_{\lambda*}(M/M_\lambda)^\sim)
\Rightarrow H^{p+q}(U^\Lambda,\cM)
\end{aligned}
$$
and we notice that $F^{pq}_2=0$ whenever $p>0$ (since
$U\cap X_{\!\lambda}$ is affine for every $\lambda\in\Lambda$)
and whenever $q>0$, since the cofiltered system
$(\Gamma(U,j_{\lambda*}(M/M_\lambda)^\sim)~|~\lambda\in\Lambda)$
has surjective transition maps. One can then argue as in the proof
of \cite[Lemma 7.3.5]{Ga-Ra} : the sheaf
$L^q:=\lim_{\lambda\in\Lambda}^q\,j_{\lambda*}(M/M_\lambda)^\sim$
is the sheafification of the presheaf : $U\mapsto H^q(U^\Lambda,\cM)$
and the latter vanishes by the foregoing. We supply an
alternative argument. Since the truly affine open subsets
form a basis of $X$, it suffices to show that $E^{0q}_2=0$
whenever $q>0$. We proceed by induction on $q$. For $q=1$,
we look at the differential $d_2^{01}:E^{0,1}_2\to E^{2,0}_2$;
by theorem \ref{th_coh-vanish} we have $E^{2,0}_2=0$, hence
$E^{0,1}_2=E^{0,1}_\infty$, and the latter vanishes by the foregoing.
Next, suppose that $q>1$, and that we have shown the vanishing
of $L^j$ for $1\leq j<q$. It follows that $E_2^{pj}=0$ whenever
$1\leq j<q$, hence $E_r^{pj}=0$ for every $r\geq 2$ and the same
values of $j$. Consequently :
$$
0=E^{0q}_\infty\simeq
\Ker(d_{q+1}^{0q}:E^{0q}_2=E^{0q}_{q+1}\to E^{q+1,0}_{q+1}).
$$
However, theorem \ref{th_coh-vanish} implies that $E^{p0}_r=0$
whenever $p>0$ and $r\geq 2$, therefore $E^{0q}_2=E^{0q}_\infty$,
{\em i.e.} $E^{0q}_2=0$, which completes the inductive step.
\end{proof}

\sset\subsubsection{}\label{subsec_closed-subs}
Let $X$ be an $\omega$-formal scheme, $\cF$ a quasi-coherent
$\cO_{\!X}$-module. For every affine open subset $U\subset X$,
we let $\mathrm{Cl}_\cF(U)$ be the set consisting of all
closed $\cO_{\!X}(U)$-submodules of $\cF(U)$. It follows
from proposition \ref{prop_quasi-coh-omega}(ii) that the rule
$$
U\mapsto\mathrm{Cl}_\cF(U)
$$
defines a presheaf on the site of all affine open subsets
of $X$. Namely, for an inclusion $U'\subset U$ of affine
open subsets, the restriction map
$\mathrm{Cl}_\cF(U)\to\mathrm{Cl}_\cF(U')$ assigns to
$N\subset\cF(U)$ the submodule $N^\sim(U')\subset\cF(U')$
(here $N^\sim$ is a quasi-coherent $\cO_{\!U}$-module).

\begin{proposition}\label{prop_closed-subs}
With the notation of \eqref{subsec_closed-subs}, the
presheaf\/ $\mathrm{Cl}_\cF$ is a sheaf on the site of
affine open subsets of $X$.
\end{proposition}
\begin{proof} Let $U\subset X$ be an affine open subset,
and $U=\bigcup_{i\in I}U_i$ a covering of $U$ by affine
open subsets $U_i\subset X$. For every $i,j\in I$ we let
$U_{ij}:=U_i\cap U_j$. Suppose there is given, for every
$i\in I$, a closed $\cO_{\!X}(U_i)$-submodule
$N_i\subset\cF(U_i)$, with the property that :
$$
N_{ij}:=N_i^\sim(U_{ij})=N^\sim_j(U_{ij})
\qquad\text{for every $i,j\in I$}
$$
(an equality of topological $\cO_{\!X}(U_{ij})$-submodules
of $\cF(U_{ij})$). Then
$(N^\sim_i)_{|U_{ij}}=(N^\sim_j)_{|U_{ij}}$, by proposition
\ref{prop_intrinsic-qcoh}, hence there exist a quasi-coherent
$\cO_{\!U}$-module $\cN$, and isomorphisms
$\cN_{|U_i}\isom N_i^\sim$, for every $i\in I$, such that
the induced $\cO_{\!U_i}$-linear maps $N^\sim_i\to\cF_{|U_i}$
assemble into a continuous $\cO_{\!U}$-linear morphism
$\phi:\cN\to\cF_{|U}$. By construction, we have a commutative
diagram of continuous maps with exact rows :
$$
\xymatrix{0 \ar[r] & \cN(U) \ar[r]^-{\rho_\cN} \ar[d] &
          \prod_{i\in I}N_i \ar[d] \ar[r] &
          \prod_{i,j\in I}N_{ij} \ar[d] \\
          0 \ar[r] & \cF(U) \ar[r]^-{\rho_\cF} &
          \prod_{i\in I}\cF(U_i) \ar[r] &
          \prod_{i,j\in I}\cF(U_{ij})
}$$
where the central vertical arrow is a closed immersion.
However, the condition of definition
\ref{def_sheaves-with-other-values}(ii) implies that both
$\rho_\cN$ and $\rho_\cF$ are admissible monomorphisms (in
the sense of \eqref{sec_Topolog-Hom-Alg}), and moreover the
image of $\rho_\cN$ is a closed submodule, since each $N_{ij}$
is a separated module. Hence also the left-most vertical arrow
is a closed immersion, and the assertion follows.
\end{proof}

\begin{definition} Let $f:X\to Y$ be a morphism of
$\omega$-formal schemes. We say that $f$ is {\em affine}
(resp. a {\em closed immersion}) if there exists a covering
$Y=\bigcup_{i\in I}U_i$ by affine open subsets, such that
for every $i\in I$, the open subset $f^{-1}U_i\subset X$
is an affine $\omega$-formal scheme (resp. is isomorphic,
as a $U_i$-scheme, to an $\omega$-formal scheme of the
form $\Spf\,A_i/J_i$, where $A_i:=\cO_{\!Y}(U_i)$ and
$J_i\subset A_i$ is a closed ideal).
\end{definition}

\begin{corollary}\label{cor_now-deserves-name}
Let $f:X\to Y$ be a morphism of $\omega$-formal schemes.
The following conditions are equivalent :
\begin{enumerate}
\item
$f$ is an affine morphism (resp. a closed immersion).
\item
For every affine open subset $U\subset Y$, the preimage
$f^{-1}U$ is an affine $\omega$-formal scheme (resp. is
isomorphic, as a $U$-scheme, to the $\omega$-formal scheme
$\Spf\,A/J$, where $A:=\cO_Y(U)$, and
$J:=\Ker(A\to\cO_{\!X}(f^{-1}U))\subset A$ is a closed ideal).
\end{enumerate}
\end{corollary}
\begin{proof} Of course, it suffices to show that
(i)$\Rightarrow$(ii). Hence, suppose that $f$ is an affine morphism,
and choose an affine open covering $Y=\bigcup_{i\in I}U_i$ such that
$f^{-1}U_i$ is an affine $\omega$-formal scheme for every $i\in I$.
We may find an affine open covering $U=\bigcup_{j\in J}V_j$ such
that, for every $j\in J$ there exists $i\in I$ with $V_j\subset
I_i$; moreover, we may assume that $J$ is countable, by proposition
\ref{prop_truly-basis}(ii.d). The assumption implies that
$f^{-1}V_j$ is affine for every $j\in J$, and then remark
\ref{rem_first-rem}(i.b) says that $A:=\cO_{\!X}(f^{-1}U)$ is an
$\omega$-admissible topological ring, so there exists a unique
morphism $g:f^{-1}U\to\Spf\,A$ such that
$g^\natural:A\to\cO_{\!X}(f^{-1}U)$ is the identity (proposition
\ref{prop_univ-prop-spf}), and $f_{|f^{-1}U}$ factors as the
composition of $g$ and the morphism $h:\Spf\,A\to U$ induced by the
natural map $B:=\cO_{\!Y}(U)\to A$. It remains to check that $g$ is
an isomorphism. To this aim, it suffices to verify that the
restriction $g_{|V_j}:f^{-1}V_j\to h^{-1}V_j$ is an isomorphism for
every $j\in J$. However, it is clear that $f_*\cO_{\!X}$ is a
quasi-coherent $\cO_{\!Y}$-module, whence a natural isomorphism of
$A$-algebras :
$$
\cO_{\!X}(f^{-1}V_j)=f_*\cO_{\!X}(V_j)\isom
A_j:=A\hat\otimes_B\cO_{\!Y}(V_j)
$$
(proposition \ref{prop_intrinsic-qcoh}) as well as an isomorphism
$f^{-1}V_j\isom \Spf\,A_j$. On the other hand,
$h^{-1}V_j=\Spf\,A_j$, and by construction $g_{|V_j}$ is the unique
morphism such that $(g_{|V_j})^\natural$ is the identity map of
$A_j$. The assertion follows.

Next, suppose that $f$ is a closed immersion, and choose
an affine open covering $Y=\bigcup_{i\in I}U_i$, and closed ideals
$J_i\subset A_i:=\cO_{\!Y}(U_i)$ such that we have isomorphisms
$f^{-1}U_i\isom\Spf\,A_i/J_i$ for every $i\in I$. Set
$U_{ij}:=U_i\cap U_j$ for every $i,j\in I$.
Clearly, $J_i^\sim(U_{ij})=J_j^\sim(U_{ij})$
for every $i,j\in I$, hence there exists a unique
closed ideal $J\subset A:=\cO_{\!Y}(U)$ such that :
\set\begin{equation}\label{eq_for-J-locally}
J^\sim(U_i)=J_i
\qquad
\text{for every $i\in I$}
\end{equation}
(proposition \ref{prop_closed-subs}). Especially, we have
$J\cO_{\!Y}(U_i)\subset J_i$, whence a unique morphism
of $U$-schemes : $g_i:f^{-1}U_i\to Z:=\Spf\,A/J$, for every
$i\in I$. The uniqueness of $g_i$ implies in particular
that $g_{i|U_{ij}}=g_{j|U_{ij}}$ for every $i,j\in I$,
whence a unique morphism $g:f^{-1}U\to Z$ of $U$-schemes.
It remains to verify that $g$ is an isomorphism, and to
this aim it suffices to check that the restriction
$g^{-1}(U_i\cap Z)\to U_i\cap Z$ is an isomorphism for
every $i\in I$. The latter assertion is clear, in view
of \eqref{eq_for-J-locally}.
\end{proof}

\begin{corollary}\label{cor_charact-truly}
Let $A$ be an $\omega$-admissible topological ring,
$U\subset X:=\Spf\,A$ an affine open subset. We have :
\begin{enumerate}
\item
The following conditions are equivalent :
\begin{enumerate}
\item
$U$ is truly affine.
\item
$U\cap\Spec\,A/I$ is quasi-compact, for every open ideal
$I\subset A$.
\end{enumerate}
\item
Especially, every quasi-compact affine open subset of $X$
is truly affine.
\end{enumerate}
\end{corollary}
\begin{proof}(i): Obviously (i.a)$\Rightarrow$(i.b). Conversely,
let $j:Y:=\Spf\,A/I\to X$ be the closed immersion, and set
$A_U:=\cO_X(U)$. On the one hand, we have $j_*\cO_Y=(A/I)^\sim$;
on the other hand, proposition \ref{prop_quasi-coh-omega}(ii)
gives an admissible short exact sequence of topological
$A_U$-modules $0\to I^\sim(U)\to A_U\to B:=(A/I)^\sim(U)\to 0$,
and notice that the topology of $B$ is discrete, since $U$ is
quasi-compact. Thus, $(Y\cap U,(\cO_Y)_{|U})=\Spf\,B=\Spec\,B$,
by corollary \ref{cor_now-deserves-name}(ii).

(ii) is an immediate consequence of (i).
\end{proof}

\begin{remark}
(i)\ \ 
Let $A$ be an $\omega$-admissible topological ring. Then
$\Spf\,A$ may well contain affine subsets that are not
truly affine. As an example, consider the one point
compactification $X:=\N\cup\{\infty\}$ of the discrete
topological space $\N$ (this is the space which induces
the discrete topology on its subset $\N$, and such that
the open neighborhoods of $\infty$ are the complements
of the finite subsets of $\N$).
We choose any field $F$, which we endow with the discrete
topology, and let $A$ be the ring of all continuous functions
$X\to F$. We endow $A$ with the discrete topology, in which
case $\Spf\,A=\Spec\,A$, and one can exhibit a natural
homeomorphism $\Spec\,A\isom X$ (exercise for the reader).
On the other hand, we have an isomorphism of topological rings :
$$
\cO_{\Spf\,A}(\N)=\lim_{b\in\N}\cO_{\Spec\,A}(\{0,\dots,b\})\simeq k^\N
$$
where $k^\N$ is endowed with the product topology. A verification
that we leave to the reader, shows that the natural injective
ring homomorphism $A\to k^\N$ induces an isomorphism of ringed
spaces
$$
\Spf\,k^\N\isom(\N,\cO_{\Spf\,A|\N})
$$
hence $\N\subset\Spf\,A$ is an affine subset. However, $\N$ is
not an affine subset of $\Spec\,A$, hence $\N$ is not a truly
affine open subset of $\Spf\,A$.

(ii)\ \
On the other hand, suppose that $A$ is $c$-admissible, in the
sense of remark \ref{rem_first-rem}(ii), and let $U\subset\Spf\,A$
be any open subset. Combining corollary \ref{cor_charact-truly}(ii)
and remark \ref{rem_first-rem}, we see that $U$ is an affine formal
scheme in the sense of \cite{EGAI-new} if and only if it is an
affine and quasi-compact $\omega$-formal scheme, if and only if
it is truly affine.
\end{remark}

We conclude this section by reviewing briefly a standard
method for getting formal schemes out of usual schemes.

\sset\subsubsection{}\label{subsec_complet-along}
Let $X$ be a scheme, $\cI\subset\cO_{\!X}$ a quasi-coherent
sheaf of ideals of finite type, and $\cF$ a quasi-coherent
$\cO_{\!X}$-module. Let also $i:X_0:=\Spec\,\cO_{\!X}/\cI\to X$
be the closed immersion. For every $n\in\N$, set
$\cF_n:=\cF/\cI^{n+1}\cF$; we endow the $i^{-1}\cO_{\!X}$-module
$i^{-1}(\cF_n)$ with its pseudo-discrete topology. Following
\cite[Ch.I, D\'ef.10.8.2]{EGAI-new}, we define the
{\em completion of $\cF$ along} $X_0$
$$
\cF^\wedge:=\lim_{n\in\N}i^{-1}\cF_n
$$
where the limit is taken in the category of sheaves of
topological abelian groups. Especially, we may consider the
completion $\cO_{\!X^\wedge}:=\cO_{\!X}^\wedge$ which is naturally
a sheaf of topological rings on $X_0$, and clearly
$\cF^\wedge$ is naturally an $\cO_{\!X^\wedge}$-module.
For every open subset $U\subset X$, set $U_0:=U\cap X_0$;
notice that, if $U\subset U'$ are two such open subsets with
$U_0=U'_0$, then the restriction map $\cF_n(U')\to\cF_n(U)$
is an isomorphism for every $n\in\N$, since the support
of $\cF_n$ lies in $X_0$. There follows a natural
identification
$$
i^{-1}\cF_n(U_0)\isom\cF_n(U)
\qquad
\text{for every $n\in\N$ and every open subset $U\subset X$}
$$
whence, by virtue of remark \ref{rem_sheaves-with-values-in-A}(iii)
a natural isomorphism of topological groups
$$
\cF^\wedge(U_0)\isom\lim_{n\in\N}\cF_n(U)
\qquad
\text{for every open subset $U\subset X$}.
$$
If moreover $U$ is quasi-compact, then $\cF_n(U)$ is a
discrete topological group for every $n\in\N$. If $U$
is affine, and $A_U:=\cO_{\!X}(U)$, $I_U:=\cI(U)$, we
have $\cF_n(U)=A_U/I_U^n\otimes_{A_U}\cF(U)$ for every
$n\in\N$, so that
$$
\cF^\wedge(U_0)=\cF(U)^\wedge
$$
where $\cF(U)^\wedge$ denotes the $I_U$-adic completion of
$\cF(U)$. Especially, the ringed space
$$
(U_0,(\cO_{\!X^\wedge})_{|U_0})
$$
is an affine $\omega$-formal scheme (see definition
\ref{def_omega-form}(iv)), isomorphic to
$\Spf\,\cO_{\!X}(U)^\wedge$. Thus
$$
X^\wedge:=(X_0,\cO_{\!X^\wedge})
$$
is an $\omega$-formal scheme called the {\em completion
of $X$ along $X_0$}, and $\cF^\wedge$ is a quasi-coherent
$\cO_{\!X^\wedge}$-module (see definition
\ref{def_quasi-coh-omega}(i)). Furthermore, the system
of projections $(i^{-1}\cF\to\cF_n)$ yields a morphism of
$i^{-1}\cO_{\!X}$-modules
$$
i^{-1}\cF\to\cF^\wedge.
$$
Especially, the map $\theta:i^{-1}\cO_{\!X}\to\cO_{\!X^\wedge}$
is a morphism of sheaves of rings, and the pair
$$
\pi_X:=(i,\theta):X^\wedge\to X
$$
is a natural morphism of ringed spaces. By inspecting the
proof of lemma \ref{lem_residue}(i), we see more precisely
that $\pi_X$ is a morphism of locally ringed spaces.

\begin{remark}\label{rem_top-on-sections}
In the situation of \eqref{subsec_complet-along}, let
$U\subset X^\wedge$ be any quasi-compact open subset. Then
the topology of $\cF(U)$ is the linear topology defined by
the system of submodules
$$
(\Gamma(U,(\cI^n\cF)^\wedge)~|~n\in\N).
$$
Indeed, the assertion is clear if $U$ is affine, since in
this case the ideal $I_U$ is finitely generated (see remark
\ref{rem_completion-of-topring}(ii,iv)). For the general
case, let $(U_\lambda~|~\lambda\in\Lambda)$ be a finite affine
covering of $U$, and endow
$M:=\prod_{\lambda\in\Lambda}\cF(U_\lambda)$ with the product
topology; according to remark \ref{rem_sheaves-with-values-in-A}(i),
the topology of $\cF(U)$ agrees with the topology induced
by $M$ via the natural injection $\cF(U)\to M$. However,
the topology of $M$ is $I_U$-adic, and we have
$\cF(U)\cap I_U^nM=\Gamma(U,(\cI^n\cF)^\wedge)$ for every
$n\in\N$, whence the contention.
\end{remark}

\sset\subsubsection{}
Keep the notation of \eqref{subsec_complet-along}, and let
$Y$ be another scheme, $\cJ\subset\cO_Y$ another quasi-coherent
sheaf of ideals, and $f:X\to Y$ a morphism of schemes such
that the associated map of sheaves $\cO_Y\to f_*\cO_{\!X}$
restricts to a morphism $\cJ\to f_*\cI$. We may then consider
the completion $Y^\wedge$ of $Y$ along $Y_0:=\Spec\cO_Y/\cJ$.
Let also $U\subset X$ and $V\subset Y$ be two affine open
subsets with $f(U)\subset V$, set
$I_U:=\cI(U)\subset A_U:=\cO_{\!X}(U)$,
$J_V:=\cJ(V)\subset B_V:=\cO_Y(V)$, and endow $A_U$ (resp.
$B_V$) with its $I_U$-adic (resp. $J_V$-adic) topology; it
follows that the resulting map $B_V\to A_U$ is continuous,
and therefore we get an induced morphism
$$
f^\wedge_{U,V}:
(U_0,(\cO_{X^\wedge})_{|U_0})\to(V_0,(\cO_{Y^\wedge})_{|V_0}).
$$
If moreover, $U'\subset U$ and $V'\subset V$ is another
pair of affine open subsets such that $f(U')\subset V'$,
it is clear that $f^\wedge_{U',V'}$ agrees with the restriction
of $f^\wedge_{U,V}$. We deduce a well defined morphism of
formal schemes
$$
f^\wedge:X^\wedge\to Y^\wedge
$$
that makes commute the diagram
$$
\xymatrix{
X^\wedge \ar[r]^-{f^\wedge} \ar[d]_{\pi_X} &
Y^\wedge \ar[d]^{\pi_Y} \\
X \ar[r]^-f & Y.
}$$
Especially, if $f$ is an open immersion and $\cI=f^{-1}\cJ$,
then clearly $f^\wedge$ is an open immersion.

\subsection{Analytically noetherian rings}
\label{sec_anal-noetherian}

\begin{definition}\label{def_anal-noetherian}
Let $A$ be an adic topological ring with a finitely
generated ideal $I$ of adic definition.

(i)\ \
We say that $A$ is {\em analytically noetherian\/} if the
following two conditions hold :
\begin{enumerate}
\alphaenu
\item
The analytic locus of $\Spec\,A$ is a noetherian scheme
(see definition \ref{def_deja-vu}).
\item
For every finitely generated $A$-module $M$, the increasing
sequence of submodules
$$
(\Ann_M(I^n)~|~n\in\N)
$$
is stationary.
\end{enumerate}

(ii)\ \
We say that $A$ is {\em universally analytically noetherian}
if, for every $n\in\N$, the polynomial $A$-algebra
$A[X_1,\dots,X_n]$, endowed with its $I$-adic topology, is
analytically noetherian.

(iii)\ \
We say that $A$ {\em satisfies the topological Artin-Rees
condition\/}, if the following holds. For every $A$-module
$M$ of finite type, and every submodule $N\subset M$, the
$I$-adic topology on $N$ agrees with the topology induced
by the $I$-adic topology of $M$.

(iv)\ \
We say that an $A$-module $M$ is {\em analytically of
finite type}, if there exists a submodule $N\subset M$
of finite type and an integer $n\in\N$ such that
$I^nM\subset N$.
\end{definition}

\begin{remark}\label{rem_anal-noetherian}
(i)\ \
It is easily seen that condition (b) of definition
\ref{def_anal-noetherian}(i) does not depend on the chosen
finitely generated ideal $I$ of adic definition of $A$.
Likewise, the class of $A$-modules of analytically finite
type is determined solely by the topology of $A$.

(ii)\ \
Obviously, every noetherian ring is analytically noetherian
when endowed with the adic topology defined by any of its
ideals.

(iii)\ \
Let $A$ be an analytically noetherian ring, and $M$ an
$A$-module of analytically finite type; then the sequence
$(\Ann_M(I^k)~|~k\in\N)$ is stationary. Indeed, pick a
finitely generated submodule $N\subset M$ such that
$I^nM\subset N$ for some $n\in\N$; by assumption, there
exists $t\in\N$ such that $\Ann_N(I^t)=\Ann_N(I^s)$ for
every $s\geq t$. It follows easily that
$\Ann_M(I^k)=\Ann_M(I^{t+n})$ for every $k\geq t+n$.
\end{remark}

\begin{lemma}\label{lem_Kato-Fuji}
Let $A$ be an adic topological ring that admits a finitely
generated ideal $I$ of adic definition. Then we have :
\begin{enumerate}
\item
The following conditions are equivalent :
\begin{enumerate}
\item
$A$ satisfies the topological Artin-Rees condition.
\item
For every $A$-module $M$ of finite type, every $n\in\N$,
and every submodule $N\subset M$ such that $I^nN=0$,
there exists $m\in\N$ such that $N\cap I^mM=0$. 
\end{enumerate}
\item
Suppose that the equivalent conditions of {\em (i)} hold
for $A$, let $M$ be any $A$-module of analytically
finite type, and $N\subset M$ any submodule. Then the
$I$-adic topology on $N$ agrees with the topology induced
by the $I$-adic topology of $M$.
\end{enumerate}
\end{lemma}
\begin{proof}(i): Clearly (a)$\Rightarrow$(b). Conversely,
let $M$ be an $A$-module of finite type, $N\subset M$
a submodule, and $n\in\N$ any integer. Set $M':=M/I^nN$
and $N':=N/I^nN$; clearly $I^nN'=0$, so (b) implies that
$N'\cap I^mM=0$ for some $m\in\N$, and the latter means
that $N\cap I^mM\subset I^nN$, whence (a).

(ii): By assumption, there exist $n\in\N$ and a submodule
$M'\subset M$ of finite type such that $I^nM\subset M'$;
set $N':=M'\cap N$. Then, for every $t\in\N$ there exists
$s\in\N$ such that
$I^sM'\cap N=I^sM'\cap N'\subset I^tN'\subset I^tN$. We
conclude that $I^{s+n}M\cap N\subset I^tN$, whence the
assertion.
\end{proof}

\begin{lemma}\label{lem_sorite-analyt}
Let $A$ be any analytically noetherian (resp. universally
analytically noetherian) topological ring, and $I\subset A$
any ideal of adic definition. We have :
\begin{enumerate}
\item
For every multiplicative subset $S\subset A$, the localization
$S^{-1}A$ is analytically noetherian (resp. universally
analytically noetherian) for its $S^{-1}I$-adic topology.
\item
Let $f:A\to B$ be a ring homomorphism, and $\cT_B$
the $IB$-adic topology on $B$. If $f$ is finite (resp.
of finite type) Then $(B,\cT_B)$ is analytically noetherian
(resp. universally analytically noetherian).
\end{enumerate}
\end{lemma}
\begin{proof}(i): Suppose that $A$ is analytically
noetherian. Clearly the analytic locus of $\Spec\,S^{-1}A$
is noetherian. Next, let $M$ be any $S^{-1}A$-module of
finite type; pick a finite system of generators
$x_\bullet:=(x_1,\dots,x_k)$ for $M$, and let $N$ be the
$A$-submodule of $M$ generated by $x_\bullet$.
Then $S^{-1}N=M$,  and $\Ann_M(S^{-1}I^n)=S^{-1}\Ann_N(I^n)$
for every $n\in\N$, whence the assertion. The assertion
for the case where $A$ is universally analytically noetherian
case follows immediately,

(ii): Again, suppose first that $A$ is analytically noetherian,
and $f$ is finite. The analytic locus of $\Spec\,B$ is finite
over the analytic locus of $\Spec\,A$, hence it is noetherian.
Next, if $M$ is a $B$-module of finite type, then it is also
an $A$-module of finite type, and $\Ann_M(I^nB)=\Ann_M(I^n)$
for every $n\in\N$, whence the assertion. Lastly, suppose that
$A$ is universally analytically noetherian and $f$ is of finite
type. For every $n\in\N$, the $A$-algebra $B[T_1,\dots,T_n]$ is
a quotient of a free polynomial $A$-algebra of finite type, so
it is analytically noetherian, by the foregoing case, whence
the assertion.
\end{proof}

\begin{proposition}\label{prop_analyt-noetherian}
Let $A$ be a ring, $I,J\subset A$ two finitely generated
ideals, and denote by $\cT_I$ (resp. $\cT_J$, resp.
$\cT_{I+J}$) the $I$-adic (resp. $J$-adic, resp. $(I+J)$-adic)
topology on $A$. We have :
\begin{enumerate}
\item
If $(A,\cT_I)$ is analytically noetherian, and $\cT_J$ is
finer than $\cT_I$, then $(A,\cT_J)$ is analytically noetherian.
\item
If $(A,\cT_I)$ and $(A,\cT_J)$ are analytically noetherian,
then the same holds for $(A,\cT_{I+J})$.
\item
If $(A,\cT_I)$ and $(A,\cT_J)$ satisfy the topological
Artin-Rees condition, the same holds for $(A,\cT_{I+J})$.
\end{enumerate}
\end{proposition}
\begin{proof}(i): After replacing $J$ by $J^n$ for a suitable
$n\in\N$, we may assume that $J\subset I$, in which case it
is clear that the analytic locus $U_I\subset X:=\Spec\,A$ of
$(A,\cT_I)$ contains the analytic locus $U_J$ of $(A,\cT_J)$.
Since by assumption $U_I$ is noetherian, the same holds for
the scheme $U_J$. Next, let $M$ be any $A$-module of finite
type, and for every $n\in\N$, denote by $\cM_n$ the
quasi-coherent $\cO_{\!X}$-module associated with $M_n:=\Ann_M(J^n)$.
Since $U_I$ is noetherian, there exists $t\in\N$ such that
$\cM_{n|U_I}=\cM_{t|U_I}$ for every $n\geq t$. Set $M':=M/M_t$,
and $M'_n:=\Ann_{M'}(J^n)$ for every $n\in\N$. It is easily
seen that
\set\begin{equation}\label{eq_Fuji}
M'_n=M_{n+t}/M_t
\qquad
\text{for every $n\in\N$}.
\end{equation}
On the other hand, since $(\cM_{n+t}/\cM_t)_{|U}=0$, we get
$$
\bigcup_{n\in\N}M'_n=N:=\bigcup_{n\in\N}\Ann_{M'}(I^n).
$$
By assumption, there exists $s\in\N$ such that
$N=\Ann_{M'}(I^s)$, and therefore $N\subset M'_s\subset N$.
Combining with \eqref{eq_Fuji} we conclude that $M_n=M_{s+t}$
for every $n\geq s+t$, whence the assertion.

(ii): Let $U_I$ and $U_J$ be as in the foregoing, and
define $U_{I+J}$ likewise as the analytic locus of
$(A,\cT_{I+J})$. Then $U_{I+J}=U_I\cup U_J$; since by
assumption $U_I$ and $U_J$ are noetherian, the same then
holds for $U_{I+J}$. Next, let $M$ be as in the foregoing;
by assumption, there exists $t\in\N$ such that
$\Ann_M(I^n)=\Ann_M(I^t)$ and $\Ann_M(J^n)=\Ann_M(J^t)$
for every $n\geq t$. It then follows easily that
$\Ann_M((I+J)^n)=\Ann_M((I+J)^{2t-1})$ for every $n\geq 2t-1$,
whence the assertion.

(iii): Let $M$ be any $A$-module of finite type, $N\subset M$
a submodule such that $(I+J)^nN=0$ for some $n\in\N$; in
light of lemma \ref{lem_Kato-Fuji}(i), it suffices to show
that there exists $m\in\N$ such that $N\cap(I+J)^mM=0$.
However, again by lemma \ref{lem_Kato-Fuji}(i), our assumption
on $(A,\cT_I)$ implies that $N\cap I^tM=0$ for some $t\in\N$;
set $M':=M/I^tM$ and let $N'\subset M'$ be the image of $N$.
Notice that $J^nN'=0$; invoking again lemma \ref{lem_Kato-Fuji}(i),
our assumption on $(A,\cT_J)$ implies that $N'\cap J^sM'=0$
for some $s\in\N$. We conclude that $N\cap(I^t+J^s)M=0$, and
then $m:=t+s-1$ will do.
\end{proof}

As a corollary, we get the following weak form of the
Artin-Rees lemma :

\begin{corollary}\label{cor_analyt-noether-TARP}
Let $(A,\cT)$ be an analytically noetherian topological ring,
$\bff:=(f_1,\dots,f_n)$ a finite sequence of elements
of $A$ that generates an ideal of adic definition. We have :
\begin{enumerate}
\item
$A$ satisfies the topological Artin-Rees condition.
\item
$A$ satisfies condition $\mathrm{(d)}_\bff$ of
\eqref{subsec_badabum}.
\end{enumerate}
\end{corollary}
\begin{proof} For every $i=1,\dots,n$, let $\cT_i$ be the
$f_iA$-adic topology on $A$; by proposition
\ref{prop_analyt-noetherian}(i), $(A,\cT_i)$ is analytically
noetherian for every $i=1,\dots,n$. Taking into account
lemma \ref{lem_sorite-analyt}(ii), we see that $A$ fulfills
the assumptions of lemma \ref{lem_Hartsho}, whence (ii).

Next, if the topological Artin-Rees condition holds for each
$(A,\cT_i)$, proposition \ref{prop_analyt-noetherian}(iii) and
a simple induction show that the same holds for $(A,\cT)$.
Thus, we may assume from start that $n=1$ and $A$ admits a
principal ideal $I=fA$ of adic definition.

Now, let $M$ be any $A$-module of finite type, $N\subset M$
a submodule, and set $M':=M/N$. Pick $t\in\N$ such that
$\Ann_{M'}(f^n)=\Ann_{M'}(f^t)$ for every $n\geq t$; it
follows easily that
$$
N\cap f^{n+t}M=f^n(N\cap f^tM)
\qquad
\text{for every $n\in\N$}
$$
whence the contention.
\end{proof}

\begin{corollary}\label{cor_analyt-noether-seseq}
Let $A$ be an analytically noetherian ring, $I\subset A$
a finitely generated ideal of adic definition, and
$$
0\to M_1\to M_2\to M_3\to 0
$$
a short exact sequence of $A$-modules, which we endow
with their $I$-adic topologies. We have :
\begin{enumerate}
\item
$M_2$ is analytically of finite type if and only if the
same holds for both $M_1$ and $M_3$.
\item
Suppose that $M_2$ is analytically of finite type.
Then the induced sequence of separated completions
\set\begin{equation}\label{eq_sep-and-complete}
0\to M_1^\wedge\to M_2^\wedge\to M_3^\wedge\to 0
\end{equation}
is short exact.
\end{enumerate}
\end{corollary}
\begin{proof}(i): Suppose that $M_2$ is analytically of
finite type. Then it is easily seen that the same holds
for $M_3$. Next, set $X:=\Spec\,A$, so that
$U:=X\setminus\Spec\,A/I$ is the analytic locus of $X$.
Let also $\cM_1$ and $\cM_2$ be the quasi-coherent
$\cO_{\!X}$-modules associated with $M_1$ and respectively
$M_2$. Then $\cM_{2|U}$ is a quasi-coherent $\cO_{\!U}$-module
of finite type, and since $U$ is noetherian, the same then
holds for its submodule $\cM_{1|U}$. Hence, we may find a
quasi-coherent $\cO_{\!X}$-submodule $\cP\subset\cM_1$ of
finite type such that $\cP_{|U}=\cM_{1|U}$, and we denote
by $P\subset M_1$ the submodule corresponding to $\cP$.
Set $M'_1:=M_1/P$ and $M'_2:=M_2/P$; by construction, the
support of $M'_1$ lies in $\Spec\,A/I$, and since $A$ is
analytically noetherian, it follows that $I^nM'_1=0$ for
some $n\in\N$. Thus, $I^nM_1\subset P$, which shows that
$M_1$ is analytically of finite type.

Lastly, suppose that $M_1$ and $M_3$ are analytically of
finite type, and pick $n\in\N$ and finitely generated
submodules $N_i\subset M_i$ such that $I^nM_i\subset N_i$
for $i=1,3$. Then we may find a submodule $P\subset M_2$
of finite type whose image in $M_3$ equals $M_3$, and
it is easily seen that $I^{2n}M_2\subset P+N_1$.

(ii): By virtue of corollary \ref{cor_analyt-noether-TARP}
and lemma \ref{lem_Kato-Fuji}(ii), the sequence
\eqref{eq_sep-and-complete} is isomorphic to the sequence
of natural maps
$$
0\to\lim_{n\in\N} M_1/(M_1\cap I^nM_2)\to
\lim_{n\in\N}M_2/I^nM_2\to\lim_{n\in\N}M_3/I^nM_3\to 0
$$
whose exactness follows easily from \cite[Lemma 3.5.3]{We}.
\end{proof}

\begin{proposition}\label{prop_like-noetherian}
Let $A$ and $I$ be as in corollary
{\em\ref{cor_analyt-noether-seseq}}, and $M$ any $A$-module
of analytically finite type. Endow $M$ with its $I$-adic
topology, and denote by $A^\wedge$ and $M^\wedge$ the separated
completions of $A$ and $M$. We have :
\begin{enumerate}
\item
The natural map $\phi_M:A^\wedge\otimes_AM\to M^\wedge$ is an
isomorphism.
\item
The completion map $A\to A^\wedge$ is flat.
\end{enumerate}
\end{proposition}
\begin{proof} We consider first the case where $M$ is of
finite type. Pick any surjective $A$-linear map $\psi:L\to M$,
from a finitely presented $A$-module $L$, and endow
$K:=\Ker\,\psi$ with its $I$-adic topology. We consider the
induced commutative diagram
$$
\xymatrix{ 0 \ar[r] & A^\wedge\otimes_AK \ar[r] \ar[d]_{\phi_K} &
A^\wedge\otimes_AL \ar[r] \ar[d]^{\phi_L} &
A^\wedge\otimes_AM \ar[r] \ar[d]^{\phi_M} & 0 \\
0 \ar[r] & K^\wedge \ar[r] & L^\wedge \ar[r] & M^\wedge \ar[r] & 0 
}$$
whose bottom row is short exact, by corollary
\ref{cor_analyt-noether-seseq}. If $L$ is free of finite
rank, $\phi_L$ is an isomorphism; then $\phi_M$ is surjective
and $\phi_K$ is injective. Since we can always find such a
surjection $\phi$ with $L$ free of finite rank, we conclude
already that $\phi_M$ is surjective for every finitely
generated $A$-module $M$.

Next, suppose that $M$ is finitely presented and $L$ is
still free of finite rank; in this case, $K$ is an $A$-module
of finite type (\cite[Lemma 2.3.18(iii)]{Ga-Ra}), so $\phi_K$
is also surjective, by the foregoing, hence $\phi_K$ is an
isomorphism and therefore the same holds for $\phi_M$.

For a general $M$ of finite type, in view of corollary
\ref{cor_analyt-noether-seseq}(i) we may find a submodule
$K'\subset K$ of finite type such that $I^nK\subset K'$.
Set $L':=L/K'$, and let $\psi':L'\to M$ be the surjective
map induced by $\psi$; then $L$' is still finitely
presented, so we may replace $L$ by $L'$ and $\psi$
by $\psi'$, and assume from start that the $I$-adic
topology is discrete on $K$, in which case the completion
map $K\to K^\wedge$ is an isomorphism, consequently
$\phi_K$ is again surjective, and we conclude as in the
foregoing that $\phi_M$ is an isomorphism, as required.

Lastly, suppose that $M$ is of analytically finite type,
and pick a submodule $N\subset M$ such that $I^n(M/N)=0$
for some $n\in\N$. We consider the analogous commutative
diagram
$$
\xymatrix{ 0 \ar[r] & A^\wedge\otimes_AN \ar[r] \ar[d]_{\phi_N} &
A^\wedge\otimes_AM \ar[r] \ar[d]^{\phi_M} &
A^\wedge\otimes_A(M/N) \ar[r] \ar[d]^{\phi_{M/N}} & 0 \\
0 \ar[r] & N^\wedge \ar[r] & M^\wedge \ar[r] & (M/N)^\wedge \ar[r] & 0 
}$$
and it is easily seen that both $\phi_{M/N}$ and the
completion map $M/N\to(M/N)^\wedge$ are isomorphisms.
The same holds for $\phi_N$, by the foregoing, so
finally $\phi_M$ is an isomorphism.

(ii): Let $N\to M$ be any injective homomorphism of
$A$-modules; we need to check that the induced map
$A^\wedge\otimes_AN\to A^\wedge\otimes_AM$ is still
injective. Since the tensor product commutes with all
colimits, we are easily reduced to the case where $M$
and $N$ are $A$-modules of finite type. In this case,
endow $M$ and $N$ with their $I$-adic topologies; by
corollary \ref{cor_analyt-noether-seseq}, the induced
map on $I$-adic completions $N^\wedge\to M^\wedge$ is
injective. Then it suffices to apply (i) to conclude.
\end{proof}

\sset\subsubsection{}\label{subsec_revised-Baire}
Let $(A,\cT)$ be a complete and separated adic
topological ring that admits a principal ideal $I=Aa$
of adic definition. Let also $M$ be an $A$-module
whose $I$-adic topology $\cT_M$ is complete and
separated. For a submodule $N\subset M$ let $N^c$
be the topological closure of $N$ in $M$, and for
every $k\in\N$ denote by $N_k\subset M$ the submodule
of all $x\in M$ such that $a^kx\in N$. Set
$$
N^s:=\bigcup_{k\in\N}N_k.
$$

\begin{lemma}\label{lem_Baire-necessities}
In the situation of \eqref{subsec_revised-Baire}, let
$N'\subset N\subset M$ be two submodules, such that
$N'$ is dense in $N$ (for the topology $\cT_M$) and
$N=N^s$. The following holds :
\begin{enumerate}
\item
The $I$-adic topology on $N$ agrees with the topology
induced by $\cT_M$.
\item
$N^c=(N^c)^s$.
\item
If\/ $N=N^c$ as well, and $N[a^{-1}]$ is an $A[a^{-1}]$-module of
finite type, then $N$ is analytically of finite type, and $N'=N$.
\end{enumerate}
\end{lemma}
\begin{proof}(i): Since $N=N^s$ we have $a^mN=a^mM\cap N$ for
every $m\in\N$, whence the assertion.

(ii): If $x\in M$ and $ax\in N^c$, then for every $n\in\N$
there exists $z_n\in M$ such that $ax-a^nz_n\in N$, and
therefore $x-a^{n-1}z_n\in N$ for every $n\geq 1$, whence
$x\in N^c$.

(ii): By assumption, there exists a finitely generated
submodule $L\subset N$ with
$$
N=L^s=\bigcup_{k\in\N}(L_k)^c
$$
(indeed, the first identity holds if we choose $L$ such that
$N[a^{-1}]=L[a^{-1}]$, and the second follows, since $N=N^c$).
Notice next that we may regard $M$ as a complete metric
space, with the metric defined by the rule
$$
d(x,y):=\left\{\begin{array}{ll}
  0     & \text{if $x=y$} \\
  2^{-b} & \text{if $x\neq y$, where $b:=\max(n\in\N~|~x-y\in a^nM)$}.
               \end{array}\right.
$$
By Baire's category theorem (\cite[Ch.IX, \S5, n.3, Th.1]{BouTG}),
it follows that $(L_k)^c$ is open in $N$ for some $k\in\N$ (for
the topology of $N$ induced by $\cT_M$), and therefore it
contains $a^mN$ for some $m\in\N$; thus
$a^{k+m}N\subset a^k(L_k)^c\subset L^c$. Combining with (i) we
deduce that the $I$-adic topology of $L^c$ agrees with the
topology induced by $\cT_M$; then by lemma
\ref{lem_tag-reinstated}(i) we conclude that $L=L^c$, so
$a^{k+m}N\subset L$, which shows that $N$ is analytically
of finite type.

Lastly, since $N'$ is dense in $N$, from (i) we see that
$N=N'+a^{m+k+1}N\subset N'+aL$. It follows that
$L=(N'\cap L)+aL$, and by Nakayama's lemma we get $L=N'\cap L$,
{\em i.e.} $L\subset N'$, so $N'$ is open in $N$ and thus
$N=N'$.
\end{proof}

\begin{theorem}
Let $(A,\cT)$ be a complete and separated adic
topological ring that admits a finitely generated ideal
$I$ of adic definition, and suppose that the analytic locus
of\/ $\Spec\,A$ is noetherian. Then $A$ is analytically
noetherian.
\end{theorem}
\begin{proof} (This is \cite[Th.5.1.2]{Fu-Ga-Ka}.) Let
$a_1,\dots,a_k$ be a finite system of elements of $A$
that generates $I$, and for every $i=1,\dots,k$ denote
by $\cT_i$ that $Aa_i$-adic topology on $A$; taking into
account proposition \ref{prop_analyt-noetherian}(ii),
it suffices to show that $(A,\cT_i)$ is analytically
noetherian for $i=1,\dots,k$. However, since the analytic
locus of $(A,\cT)$ is noetherian, the same holds for the
analytic locus of $(A,\cT_i)$, and the topology $\cT_i$
is also complete and separated for every $i=1,\dots,k$,
by virtue of lemma \ref{lem_fontaine}. We are then
reduced to the case where $I=Aa$ is a principal ideal.

Now, let $M$ be an $A$-module of finite type; we need to
show that the sequence of submodules $(\Ann_M(a^n)~|~n\in\N)$
is stationary. Let us write $M=L/Q$ for a free $A$-module
$L$ of finite rank, and a submodule $Q\subset L$, and
set $N:=Q^s$ (notation of \eqref{subsec_revised-Baire}).
Notice that
$$
N/Q=\bigcup_{n\in\N}\Ann_M(a^n)
$$
so it suffices to show that $a^kN\subset Q$ for some $k\in\N$.
Now, obviously $N^s=N$, so $N^c=(N^c)^s$ as well, by lemma
\ref{lem_Baire-necessities}(ii). Moreover, since $A[a^{-1}]$
is noetherian, $N^c[a^{-1}]$ is an $A[a^{-1}]$-module of finite
type, and consequently $N^c$ is analytically of finite type,
and $N=N^c$, by lemma \ref{lem_Baire-necessities}(ii).
It follows that $N/Q$ is also analytically of finite type;
so, pick a submodule $N'\subset N/Q$ of finite type and an
integer $k\in\N$ such that $a^k(N/Q)\subset N'$; since $N'$
is finitely generated, there exists as well $n\in\N$ such
that $a^nN'=0$, so finally $a^{k+n}(N/Q)=0$, as required.
\end{proof}

Next, we globalize definition \ref{def_anal-noetherian}
to schemes as follows.

\begin{definition}\label{def_analytic-schemes}
Let $A$ be an adic topological ring that has a finitely
generated ideal $I$ of adic definition, $X$ an $A$-scheme,
and $\cF$ a quasi-coherent $\cO_{\!X}$-module.

(i)\ \
We say that $X$ is {\em locally analytically noetherian\/}
if every point of $X$ admits an affine open neighborhood $U$
such that the following holds. Endow $A_U:=\cO_{\!X}(U)$ with
the $I$-adic topology $\cT_U$; then $(A_U,\cT_U)$ is analytically
noetherian.

(ii)\ \
We say that $X$ is {\em analytically noetherian\/} if it
is quasi-compact, quasi-separated and locally analytically
noetherian.

(iii)\ \
We say that $\cF$ is {\em analytically of finite type\/} if
every point of $X$ admits an affine open neighborhood $U$ such
that $\cF(U)$ is an $(A_U,\cT_U)$-module of analytically finite
type.
\end{definition}

\begin{remark}\label{rem_an-noether-sch}
If $A$ is a universally analytically noetherian ring,
then every $A$-scheme of finite type is analytically
noetherian, by virtue of lemma \ref{lem_sorite-analyt}.
\end{remark}

\begin{lemma}\label{lem_global-approx}
Let $A$ and $I$ be as in definition {\em\ref{def_analytic-schemes}}.
Let also $X$ be a quasi-compact and quasi-separated $A$-scheme,
and $\cF$ a quasi-coherent $\cO_{\!X}$-module. We have :
\begin{enumerate}
\item
If $\cF$ is analytically of finite type, there exist an
integer $n\in\N$ and a quasi-coherent $\cO_{\!X}$-submodule
of finite type $\cG\subset\cF$ such that $I^n\cF\subset\cG$.
\item
If $\cF$ is analytically of finite type and $X$ is analytically
noetherian, there exists $p\in\N$ such that
$\Ann_\cF(I^q)=\Ann_\cF(I^p)$ for every integer $q\geq p$.
\item
If $\cF$ is of finite type and $X$ is analytically
noetherian, there exist a finitely presented
$\cO_{\!X}$-module $\cG$, an epimorphism $\phi:\cG\to\cF$
of $\cO_{\!X}$-modules, and an integer $n\in\N$ such
that $I^n\Ker\,\phi=0$.
\end{enumerate}
\end{lemma}
\begin{proof}(i): By assumption, we may find a finite affine open
covering $(U_\lambda~|~\lambda\in\Lambda)$ of $X$, and for every
$\lambda\in\Lambda$ a quasi-coherent $\cO_{\!U_\lambda}$-submodule
$\cH_\lambda$ of finite type of $\cF_\lambda:=\cF_{|U_\lambda}$ and an
integer $n_\lambda\in\N$ such that $I^{n_\lambda}\cF_\lambda\subset\cH_\lambda$.
We may also find, for every $\lambda\in\Lambda$, a finitely
presented quasi-coherent $\cO_{\!U_\lambda}$-module $\cH'_\lambda$ with
an $\cO_{\!U_\lambda}$-linear morphism
$\phi_\lambda:\cH'_\lambda\to\cF_\lambda$ whose image is $\cH_\lambda$.
According to lemma \ref{lem_extend-cohs}(i) there exist a
finitely presented quasi-coherent $\cO_{\!X}$-module $\cG_\lambda$
with $\cG_{\lambda|U_\lambda}=\cH'_\lambda$ and a morphism
$\psi_\lambda:\cG_\lambda\to\cF_\lambda$ of $\cO_{\!X}$-modules
such that $\psi_{\lambda|U_\lambda}=\phi_\lambda$. Then we may take
$\cG:=\sum_{\lambda\in\Lambda}\Img\,\psi_\lambda$.

(ii): Suppose that for every $\lambda\in\Lambda$ there exists
$p_\lambda\in\N$ such that
$\Ann_{\cF_\lambda}(I^q)=\Ann_{\cF_\lambda}(I^{p_\lambda})$ for every
$q\geq p_\lambda$. Then $p:=\max(p_\lambda~|~\lambda\in\Lambda)$
fulfills the stated condition. Thus, we may suppose
that $X=\Spec\,A$ for an analytically noetherian ring
$A$, in which case the assertion follows from remark
\ref{rem_anal-noetherian}(iii).

(iii): By assumption, we may find a finite affine covering
$(U_\lambda~|~\lambda\in\Lambda)$ of $X$ such that
$A_\lambda:=\cO_{\!X}(U_\lambda)$ is analytically noetherian
relative to its $IA_\lambda$-adic topology, for every
$\lambda\in\Lambda$. Pick a finitely presented quasi-coherent
$\cO_{\!U_\lambda}$-module $\cH_\lambda$ and an epimorphism
$\psi_\lambda:\cH_\lambda\to\cF_{|U_\lambda}$ for every
$\lambda\in\Lambda$. By lemma \ref{lem_extend-cohs}(i)
there exist a finitely presented quasi-coherent
$\cO_{\!X}$-module $\cH'_\lambda$ with
$\cH'_{\lambda|U_\lambda}=\cH_\lambda$ and a morphism of
$\cO_{\!X}$-modules $\psi'_\lambda:\cH'_\lambda\to\cF$ such
that $\psi'_{\lambda|U_\lambda}=\psi_\lambda$, for every
$\lambda\in\Lambda$. Set
$\cH:=\bigoplus_{\lambda\in\Lambda}\cH'_\lambda$; the sum of
the morphisms $\psi'_\lambda$ is an epimorphism $\psi:\cH\to\cF$,
and $\cH$ is clearly finitely presented. For every
$\lambda\in\Lambda$, let also $V_\lambda\subset U_\lambda$
be the analytic locus; then $V:=\bigcup_{\lambda\in\Lambda}V_\lambda$
is locally noetherian, quasi-compact and quasi-separated,
and therefore $\cK:=\Ker\,\psi$ restricts to a finitely
presented quasi-coherent $\cO_V$-module $\cK_{|V}$. Then,
invoking again lemma \ref{lem_extend-cohs}(i) we find a
quasi-coherent finitely presented $\cO_{\!X}$-module
$\cK'$ and an $\cO_{\!X}$-linear morphism $\cK'\to\cK$
that restricts to an isomorphism on $V$. Let $\cK''\subset\cK$
be the image of $\cK'$, and set $\cG:=\cH/\cK''$. Clearly
$\psi$ factors through an epimorphism $\phi:\cG\to\cF$,
and the support of $\Ker\,\phi$ lies in $X\setminus V$.
Lastly, since $A_\lambda$ is analytically noetherian, we
may find an integer $n_\lambda\in\N$ such that
$I^{n_\lambda}\cdot(\Ker\,\phi)_{|V_\lambda}=0$, for every
$\lambda\in\Lambda$. The assertion then holds with
$n:=\max(n_\lambda~|~\lambda\in\Lambda)$.
\end{proof}

\begin{proposition}\label{prop_glob-analyt-noeth}
Let $A$ be an adic topological ring with a finitely generated
ideal $I$ of adic definition, $B$ an $A$-algebra, that we endow
with its $IB$-adic topology $\cT_B$. Let also $F$ be any $B$-module,
and denote by $\cF$ the quasi-coherent $\cO_{\Spec\,B}$-module
arising from $F$. We have :
\begin{enumerate}
\item
The following conditions are equivalent :
\begin{enumerate}
\alphaenu
\item
$(B,\cT_B)$ is analytically noetherian.
\item
$\Spec\,B$ is an analytically noetherian $A$-scheme.
\end{enumerate}
\item
The following conditions are equivalent :
\begin{enumerate}
\item
$F$ is a $B$-module of analytically finite type.
\item
$\cF$ is an $\cO_{\Spec\,B}$-module of analytically finite type.
\end{enumerate}
\end{enumerate}
\end{proposition}
\begin{proof}(i): Clearly (a)$\Rightarrow$(b). Hence, suppose
that $X:=\Spec\,B$ is analytically noetherian, and let
$(U_i~|~i=1,\dots,k)$ be a finite affine covering of $X$
such that $B_i:=\cO_{\!X}(U_i)$ is analytically noetherian
for its $IB_i$-adic topology for every $i=1,\dots,k$. By
virtue of lemma \ref{lem_sorite-analyt}, we may assume
that for every $i=1,\dots,k$ there exists $b_i\in B$ such
that $B_i=B[b_i^{-1}]$. Clearly the analytic locus $U$ of
$X$ is the union of the analytic loci of
$\Spec\,B_1,\dots,\Spec\,B_k$, so $U$ is noetherian.
Next, let $M$ be any $B$-module of finite type, and
for every $i=1,\dots,k$, set $M_i:=B_i\otimes_BM$;
choose $n\in\N$ such that $\Ann_{M_i}(I^rB_i)=\Ann_{M_i}(I^nB_i)$
for every $i=1,\dots,n$ and every $r\geq n$. It follows
that $\Ann_M(I^r)=\Ann_M(I^n)$ for every $r\geq n$,
whence the contention.

(ii) follows easily from lemma \ref{lem_global-approx}(i).
\end{proof}

\sset\subsubsection{}\label{subsec_complet-univ-an-noeth-sch}
Let $A$ be a universally analytically noetherian ring,
$I\subset A$ a finitely generated ideal of adic definition,
$X$ a quasi-separated $A$-scheme locally of finite type,
and $\cF$ a quasi-coherent $\cO_{\!X}$-module. According to
\eqref{subsec_complet-along} we may consider the
completions $X^\wedge$ and $\cF^\wedge$ of $X$ and $\cF$
along the closed subscheme
$X_0:=\Spec\,\cO_{\!X}/I\cO_{\!X}=\Spec\,A/I\times_{\Spec\,A}X$,
which are respectively a formal scheme endowed with a natural
morphism of locally ringed spaces
$$
\pi:X^\wedge\to X
$$
and a quasi-coherent $\cO_{\!X^\wedge}$-module with a natural
map of $\cO_{\!X^\wedge}$-modules
\set\begin{equation}\label{eq_complete-module}
\pi^*\cF\to\cF^\wedge.
\end{equation}

\begin{proposition}\label{prop_up-to-completion}
With the notation of \eqref{subsec_complet-univ-an-noeth-sch},
the following holds :
\begin{enumerate}
\item
If $\cF$ is analytically of finite type, the map
\eqref{eq_complete-module} is an isomorphism.
\item
We have $H^i(U,\pi^*\cF)=0$ for every affine quasi-compact
open subset $U\subset X^\wedge$ and every $i>0$.
\item
If $X$ is affine, the natural map
$$
\cO_{\!X^\wedge}(X^\wedge)\otimes_{\cO_{\!X}(X)}H^0(X,\cF)\to
H^0(X^\wedge,\pi^*\cF)
$$
is an isomorphism.
\item
Let $\fU:=(U_i~|~i\in I)$ be an open covering of $X^\wedge$,
and suppose that $U_i\cap U_j$ is affine for every $i,j\in I$.
Then there is a natural isomorphism (notation of
\eqref{subsec_altern-pseudo-Leray}) :
$$
H^\bullet_\alt(\fU,\pi^*\cF)\isom H^\bullet(X^\wedge,\pi^*\cF).
$$
\end{enumerate}
\end{proposition}
\begin{proof}(i): The assertion is local on $X$, hence
we may assume that $X$ is local, say $X=\Spec\,B$, and
$\cF$ is the quasi-coherent $\cO_{\!X}$-module attached
to the module $F:=\cF(X)$. Let $T$ be any final object
in the category of topological spaces ({\em i.e.} $T$
is the unique topology over a set that has only one
point), and endow $T$ with a structure sheaf $\cO_T$
by declaring that $\cO_T(T):=B$ and $\cO_T(\emptyset):=0$.
The module $F$ defines an $\cO_T$-module $F_T$ by the
rule : $F_T(T):=F$ and $F_T(\emptyset):=0$. The system
of restriction maps $B\to\cO_{\!X}(U)$ for $U$ ranging
over the open subsets of $X$, may be viewed as a morphism
of ringed spaces
$$
\eps:X\to T
$$
and with this notation we then have $\cF=\eps^*F_T$.
Set as well $\eps^\wedge:=\eps\circ\pi$; there follows
a natural isomorphism of $\cO_{\!X^\wedge}$-modules
$$
\pi^*\cF\isom\eps^{\wedge*}F_T.
$$
In other words, $\pi^*\cF$ is the sheaf associated to the
presheaf given by the rule :
$$
U_0\mapsto\cO_{\!X^\wedge}(U_0)\otimes_BF
\qquad
\text{for every open subset $U_0\subset X^\wedge=X_0$}
$$
and the discussion of \eqref{subsec_complet-along} shows
that \eqref{eq_complete-module} is the morphism associated
to the morphism of presheaves given, on every affine
open subset $U_0\subset X_0$, by the natural map
\set\begin{equation}\label{eq_nat-to-completion}
\cO_{\!X}(U)^\wedge\otimes_BF\to\cF(U)^\wedge
\end{equation}
where $U\subset X$ is any affine open subset such that
$U_0=U\cap X_0$, and where $\cO_{\!X}(U)^\wedge$ and
$\cF(U)^\wedge$ denote the $I$-adic completions of
$\cO_{\!X}(U)$ and respectively $\cF(U)$. However, since
$\cF$ is quasi-coherent, we have
$\cF(U)=\cO_{\!X}(U)\otimes_BF$; taking into account
remark \ref{rem_an-noether-sch} and propositions 
\ref{prop_like-noetherian} and \ref{prop_glob-analyt-noeth}(ii)
we deduce that \eqref{eq_nat-to-completion} is an isomorphism
for every such $U_0$. Since the affine open subsets are
a basis of the topology of $X^\wedge$ that is closed under
finite intersections, the assertion follows.

(ii): Since $U$ is quasi-compact, $\pi(U)$ is contained
in a quasi-compact open subset $V\subset X$; we may then
replace $X$ by $V$, and assume from start that $X$ is
quasi-compact and quasi-separated. In this case, $\cF$
is the colimit of the filtered family
$(\cF_\lambda~|~\lambda\in\Lambda)$ of its quasi-coherent
$\cO_{\!X}$-submodules of finite type (proposition
\ref{prop_fp-approx}), hence $\pi^*\cF$ is the colimit
of the system $(\pi^*\cF_\lambda~|~\lambda\in\Lambda)$.
Since $U$ is a spectral topological space, proposition
\ref{prop_dir-im-and-colim} then reduces to showing that
$H^i(U,\pi^*\cF_\lambda)=0$ for every $\lambda\in\Lambda$.
Hence, we may assume from start that $\cF$ is of finite
type, in which case $\pi^*\cF$ is a quasi-coherent
$\cO_{\!X^\wedge}$-module, by (i), and then the assertion
follows from theorem \ref{th_coh-vanish}(ii).

(iii): Arguing as in the proof of (ii), we reduce to the
case where $\cF$ is analytically of finite type; then
$\cF(X)$ is an $\cO_{\!X}(X)$-module of analytically
finite type (proposition \ref{prop_glob-analyt-noeth}(ii)),
and the discussion of \eqref{subsec_complet-along} says
that $\cO_{\!X^\wedge}(X^\wedge)$ (resp. $\cF^\wedge(X^\wedge)$) is
the $I$-adic completion of $\cO_{\!X}(X)$ (resp. of $\cF(X)$).
Then the assertion follows from (i), lemma \ref{lem_left-adj}(i)
and proposition \ref{prop_like-noetherian}(i).

(iv) follows from (ii), together corollary \ref{cor_Leray}(ii)
and \cite[Ch.I, Prop.10.7.2]{EGAI-new}.
\end{proof}

\begin{theorem}\label{th_analyt-proper-finiteness}
Let $A$ be a universally analytically noetherian ring,
$X$ a proper and finitely presented $A$-scheme, and $\cF$
a quasi-coherent $\cO_{\!X}$-module of analytically finite
type. Then the $A$-module $H^i(X,\cF)$ is analytically
of finite type for every $i\in\N$.
\end{theorem}
\begin{proof} Let $I\subset A$ be any ideal of adic
definition. The first observation is the following :

\begin{claim}\label{cl_reduce-to-fpres}
Let $i\in\N$ be any integer, and suppose that $H^i(X,\cF')$
is analytically of finite type for every finitely presented
$\cO_{\!X}$-module $\cF'$. Then $H^i(X,\cF)$ is analytically
of finite type for every $\cO_{\!X}$-module $\cF$ of
analytically finite type.
\end{claim}
\begin{pfclaim} Let $\cG\subset\cF$ be an $\cO_{\!X}$-submodule
as in lemma \ref{lem_global-approx}(i), so that $I^n\cF\subset\cG$
for some $n\in\N$, and denote by $j:\cG\to\cF$ the inclusion
morphism. Then, for every $a\in I^n$ the endomorphism
$a\cdot\one_\cF$ factors thorugh an $\cO_{\!X}$-linear morphism
$\phi_a:\cF\to\cG$. Set $H^i:=H^i(X,\cF)$; we deduce that
$a\cdot\one_{H^i}=H^i(X,j)\circ H^i(X,\phi_a)$.
Especially, $I^nH^i$ is contained in the image of $H^i(X,j)$.
Hence, if $H^i(X,\cG)$ is analytically of finite type, the
same follows for $H^i$. We may then replace $\cF$ by $\cG$,
and assume from start that $\cF$ is quasi-coherent of finite
type. Next, by lemma \ref{lem_global-approx}(iii) and remark
\ref{rem_an-noether-sch} we may find a short exact sequence of
quasi-coherent $\cO_{\!X}$-modules $0\to\cK\to\cF'\to\cF\to 0$
such that $\cF'$ is finitely presented and $I^n\cK=0$ for some
$n\in\N$. It follows easily that $I^n$ annihilates the cokernel
of the induced map $H^i(X,\cF')\to H^i$, whence the claim.
\end{pfclaim}

Thus, suppose henceforth that $\cF$ is finitely presented,
and set $S:=\Spec\,A$; according to
\cite[Ch.IV, Prop.8.9.1(i,ii)]{EGAIV-3} we may assume
that there exist a noetherian subring $A_0\subset A$, a
finitely presented morphism of schemes $X_0\to S_0:=\Spec\,A_0$
and a coherent $\cO_{\!X_0}$-module $\cF_0$ such that
$X=S\times_{S_0}X_0$ and $\cF=\pi^*\cF_0$, where $\pi:X\to X_0$
is the projection. By \cite[Ch.IV, Th.8.10.5]{EGAIV-3}, we may
moreover assume that $X_0$ is a proper $A_0$-scheme.

\begin{claim}\label{eq_local-Tors}
The $\cO_{\!X}$-module $\cTor_q^{S_0}(\cO_{\!S},\cF_0)$
is analytically of finite type for every $q\in\N$.
\end{claim}
\begin{pfclaim} The assertion is local on $X_0$, so we
may assume that $X_0$ is and affine and finitely presented
$S_0$-scheme; then we may find an integer $n\in\N$ and a
closed immersion $i:X_0\to\A^n_{S_0}$. Set
$i_S:=S\times_{S_0}i:X\to\A^n_S$; according to
\cite[Ch.III, Prop.6.5.11]{EGAIII-2} there is a natural
isomorphism of $\cO_{\!X}$-modules :
$$
i_{S*}\cTor_q^{S_0}(\cO_{\!S},\cF_0)\isom
\cTor_q^{S_0}(\cO_{\!S},i_*\cF_0)
$$
so we may replace $X_0$ by $\A^n_{S_0}$ and $\cF_0$ by
$i_*\cF_0$, and assume from start that $X_0$ is a flat
and affine $S_0$-scheme. In this situation, there exists
as well a natural isomorphism of $\cO_{\!X}$-modules :
$$
\cTor_q^{S_0}(\cO_{\!S},\cF_0)\isom
\cT_q:=\cTor_q^{X_0}(\cO_{\!X},\cF_0)
\qquad
\text{for every $n\in\N$}.
$$
Indeed, for every $V_0$, $V$ and $U$ as in the foregoing
condition (a) we have an isomorphism
$$
\gamma_{V,U}:\Tor^{\cO_{\!S_0}(V_0)}_q(\cO_{\!S}(V),\cF_0(U))\isom
\Tor^{\cO_{\!X_0}(U)}_q(\cO_{\!X}(V\times_{S_0}U),\cF_0(U))
$$
of $\cO_{\!X}(V\times_{S_0}U)$-modules, such that for every
inclusion of open subsets $V'_0\subset V_0$,
$V'\subset V\cap(S\times_{S_0}V'_0)$ and
$U'\subset U\cap(S\times_{S_0}V'_0)$ the resulting diagram
commutes :
$$
\xymatrix{ \Tor^{\cO_{\!S_0}(V_0)}_q(\cO_{\!S}(V),\cF_0(U))
\ar[rr]^-{\gamma_{V,U}} \ar[d] & &
\Tor^{\cO_{\!X_0}(U)}_q(\cO_{\!X}(V\times_{S_0}U),\cF_0(U))  \ar[d] \\
\Tor^{\cO_{\!S_0}(V'_0)}_q(\cO_{\!S}(V'),\cF_0(U')) \ar[rr]^-{\gamma_{V',U'}} & &
\Tor^{\cO_{\!X_0}(U')}_q(\cO_{\!X}(V'\times_{S_0}U'),\cF_0(U'))
}$$
(\cite[Prop.3.2.9]{We}) whence the assertion. Thus, it
suffices to check that $\cT_q$ is analytically of finite
type; however, since $X_0$ is affine, $\cF_0$ admits
a resolution $\cP_\bullet\isom\cF[0]$ consisting of free
$\cO_{\!X_0}$-modules of finite type, and $\cT_q$ is
isomorphic to $H_q(\cO_{\!X}\otimes_{\cO_{\!X_0}}\cP_\bullet)$.
The claim then follows easily from corollary
\ref{cor_analyt-noether-seseq}(i).
\end{pfclaim}

Let $\fU:=(U_i~|~i=1,\dots,n)$ be a finite affine
covering of $X_0$; we have natural isomorphisms in
$\sD(A_0\Mod)$ and respectively $\sD(A\Mod)$
$$
R\Gamma(X_0,\cF_0)\isom\bar C{}^\bullet_\mathrm{alt}(\fU,\cF_0)
\qquad
R\Gamma(X,\cF)\isom A\otimes_{A_0}\bar C{}^\bullet_\mathrm{alt}(\fU,\cF_0)
$$
where $\bar C{}^\bullet_\mathrm{alt}(\fU,-)$ denotes the truncated
alternating \v{C}ech complex associated with the covering $\fU$
(theorem \ref{th_Cech-resolve}(ii)).
Let us choose a Cartan-Eilenberg projective resolution
$P^{\bullet\bullet}\isom\bar C{}^\bullet_\mathrm{alt}(\fU,\cF_0)[0]$
such that $P^{pq}=0$ whenever $q>0$ (see \cite[Lemma 5.7.2]{We});
the double complex $A\otimes_{A_0}P^{\bullet\bullet}$ gives
rise to a spectral sequence
$$
E^{p,-q}_2:=\Tor_q^{A_0}(A,H^p\bar C{}^\bullet_\mathrm{alt}(\fU,\cF_0))
\Rightarrow
M_{p-q}:=H_{p-q}(A\derotimes_{A_0}\bar C{}^\bullet_\mathrm{alt}(\fU,\cF_0))
$$
as well as a spectral sequence
$$
F^{p,-q}_1:=\Tor_q^{A_0}(A,\bar C{}^p_\mathrm{alt}(\fU,\cF_0))
\Rightarrow M_{p-q}
$$
whose differentials $d_1^{pq}:F^{pq}_1\to F_1^{p+1,q}$ are induced
by those of $\bar C{}^\bullet_\mathrm{alt}(\fU,\cF_0)$, whence natural
isomorphisms for every $p,q\in\N$ :
\set\begin{equation}\label{eq_sprouts}
E^{p,-q}_2\isom\Tor_q^{A_0}(A,H^p(X_0,\cF_0))
\qquad
F^{p,-q}_2\isom
H^p(X,\cTor^{S_0}_q(\cO_{\!S},\cF_0))
\end{equation}
and especially, we have natural isomorphisms
\set\begin{equation}\label{eq_llama}
F^{i,0}_2\isom H^i
\qquad
\text{for every $i\in\N$}.
\end{equation}

\begin{claim}\label{cl_abut-an-ftype}
$M_{p-q}$ is analytically of finite type, for every $p,q\in\N$.
\end{claim}
\begin{pfclaim} The $A_0$-module $H^p(X_0,\cF_0)$ is
finitely generated for every $p\in\N$, because $X_0$ is a
proper $A_0$-scheme (\cite[Ch.III, Th.3.2.1]{EGAIII});
since $A_0$ is noetherian, we may then find a resolution
$P_\bullet$ of $H^p(X_0,\cF_0)$ consisting of free $A_0$-modules
of finite type, and there are natural isomorphisms
$E^{p,-q}_2\isom H_q(A\otimes_{A_0}P_\bullet)$; taking into
account corollary \ref{cor_analyt-noether-seseq}(i),
we conclude that $E^{p,-q}_2$ is an $A$-module of analytically
finite type, for every $p,q\in\N$. Then, again by corollary
\ref{cor_analyt-noether-seseq}(i), a simple induction on
$r\in\N$ shows that $E^{p,-q}_r$ is analytically of finite
type, for every $p,q,r\in\N$. We deduce that the $A$-module
$M_{p-q}$ admits a finite filtration whose subquotients are
analytically of finite type, and to conclude it suffices
to invoke again corollary \ref{cor_analyt-noether-seseq}(i).
\end{pfclaim}

Now the theorem is a special case of the following more general

\begin{claim} $F^{p,-q}_{r+2}$ is an $A$-module of analytically
finite type, for every $p,q,r\in\N$.
\end{claim}
\begin{pfclaim}[] We argue by descending induction on $p$,
and notice that $F^{p,-q}_1=0$ for every $p\geq n-1$, so the
assertion is clear for $p\geq n-1$. Next, let $t\leq n-1$,
and suppose that the claim has already been shown for every
$p,q,r\in\N$ with $r\geq 2$ and $p\geq t$.
First we prove, by descending induction on $r\geq 2$, that
$F^{t-1,0}_r$ is analytically of finite type. Indeed, there
exists $s\in\N$ large enough, such that $F^{t-1,0}_r$ is a
subquotient of $M_{t-1}$ for every $r\geq s$, whence the
assertion for every $r\geq s$, by virtue of claim
\ref{cl_abut-an-ftype} and corollary
\ref{cor_analyt-noether-seseq}(i). Thus, suppose that
$2<r\leq s$, and we know already that $F^{t-1,0}_r$
is analytically of finite type; we have an exact
sequence
$$
0\to F^{t-1,0}_r\to F^{t-1,0}_{r-1}\to F^{t-2+r,2-r}_{r-1}
$$
and by inductive assumption we also know that $F^{t-2+r,2-r}_{r-1}$
is analytically of finite type. From corollary
\ref{cor_analyt-noether-seseq}(i), it follows that
the same holds for $F^{t-1,0}_{r-1}$, as required.
Especially, in light of \eqref{eq_llama}, the foregoing
implies that $H^{t-1}(X,\cF)$ is analytically of finite type,
for every finitely presented $\cO_{\!X}$-module $\cF$, hence
the same holds more generally whenever $\cF$ is analytically
of finite type, by virtue of claim \ref{cl_reduce-to-fpres}.
Taking into account claim \ref{eq_local-Tors} and
\eqref{eq_sprouts}, we deduce that $F^{t-1,-q}_2$ is
analytically of finite type for every $q\in\N$, and
finally, the same follows for $F^{t-1,-q}_r$, whenever
$q,r\in\N$ and $r\geq 2$, after invoking once more
corollary \ref{cor_analyt-noether-seseq}(i).
\end{pfclaim}
\end{proof}

\sset\subsubsection{}\label{subsec_filtration-on-coh}
In the situation of theorem \ref{th_analyt-proper-finiteness},
let $I\subset A$ be any finitely generated ideal of adic
definition; we notice that $I^n\cF$ is a quasi-coherent
$\cO_{\!X}$-module of analytically finite type for every
$\in\N$, and we define a descending filtration
$\Fil_I^\bullet H^i$ on $H^i:=H^i(X,\cF)$ by setting
$$
\Fil^n_IH^i:=\Img(\psi_{i,n}:H^i(X,I^n\cF)\to H^i)
\qquad
\text{for every $i,n\in\N$}.
$$

\begin{corollary}\label{cor_analyt-finiteness}
With the notation of \eqref{subsec_filtration-on-coh}, the
following holds for every $i\in\N$ :
\begin{enumerate}
\item
The linear topology on $H^i$ defined by the descending
filtration $\Fil_I^\bullet H^i$ agrees with the $I$-adic
topology.
\item
The system $(\Ker\,\psi_{i,n}~|~n\in\N)$ is essentially zero.
\end{enumerate}
\end{corollary}
\begin{proof} Let us show first that (i)$\Rightarrow$(ii).
To ease notation, set $K^i_n:=\Ker\,\psi_{i,n}$ for every
$i,n\in\N$. For every $a\in I^n$, the endomorphism
$a\cdot\one_{\cF}$ is the composition of the inclusion map
$j_n:I^n\cF\to\cF$ and a morphism $\phi_{a,n}:\cF\to I^n\cF$ of
$\cO_{\!X}$-modules, and we have as well
$$
\phi_{a,n}\circ j_n=a\cdot\one_{I^n\cF}
$$
whence $H^i(X,\phi_{a,n})\circ\psi_{i,n}=a\cdot\one_{H^i(X,I^n\cF)}$.
We deduce immediately that
\set\begin{equation}\label{eq_godspeed}
I^nK_n^i=0
\qquad
\text{for every $p,n\in\N$}.
\end{equation}
Next, notice that
$$
\Img(K^i_p\to K^i_n)=K^i_n\cap\Img(H^i(X,I^p\cF)\to H^i(X,I^n\cF))
\qquad
\text{for all integers $p\geq n\geq 0$}.
$$
Suppose now that (i) holds for every $\cO_{\!X}$-module $\cF$
of analytically finite type; especially, it applies to
$I^n\cF$, and therefore for every $k\in\N$ we may find an
integer $p\geq n$ such that
$\Img(H^i(X,I^p\cF)\to H^i(X,I^n\cF))\subset I^kH^i(X,I^n\cF)$,
whence
$$
\Img(K^i_p\to K^i_n)\subset K^i_n\cap I^kH^i(X,I^n\cF).
$$
By theorem \ref{th_analyt-proper-finiteness}, we know
that $H^i(X,I^n\cF)$ is analytically of finite type;
by corollary \ref{cor_analyt-noether-TARP} and lemma
\ref{lem_Kato-Fuji}(ii), it follows that there exists
$k\in\N$ large enough, such that
$K^i_n\cap I^kH^i(X,I^n\cF)\subset I^nK^i_n$, and then
the assertion follows from \eqref{eq_godspeed}.

(i): Let $a_1,\dots,a_k$ be a finite system of generators
for $I$; we argue by induction on $k$. If $k=1$, let
$p\in\N$ such that $\Ann_\cF(a^q_1)=\Ann_\cF(a^p_1)$ for every
$q\geq p$ (lemma \ref{lem_global-approx}(ii)); there
follows a commutative diagram of $\cO_{\!X}$-modules
$$
\xymatrix{ a_1^p\cF \ar[rr] \ar[d] & & a_1^q\cF \ar[d] \\
\cF \ar[rr]^-{a_1^{q-p}\cdot\one_\cF} & & \cF
}$$
whose top horizontal arrow is an isomorphism, and
whose vertical arrows are the inclusion maps. We deduce
easily that $\Fil_I^qH^i=a_1^{q-p}\cdot\Fil_I^pH^i$ for every
$q\geq p$, whence
$$
a_1^qH^i\subset\Fil_I^qH^i\subset a_1^{q-p}H^i
$$
as required. Next, let $k>1$, and set $J:=Aa_1+\dots+Aa_{k-1}$
and $L:=Aa_k$. Denote by $\cT_I$, $\cT_J$ and $\cT_L$ respectively
the $I$-adic, $J$-adic and $L$-adic topologies on $A$.

Fix $p\in\N$; it is easily seen that $I^pH^i\subset\Fil^pH^i$,
so it remains only to check that there exists $n\in\N$ such
that $\Fil^nH^i\subset I^pH^i$. To this aim, notice first that
$X$ is also an analytically noetherian $(A,\cT_L)$-scheme
(proposition \ref{prop_analyt-noetherian}(i)), and $\cF$
is also analytically of finite type relative to the topology
$\cT_L$; since $L$ is principal, it follows that we may find
$t\in\N$ such that
\set\begin{equation}\label{eq_almost-over}
\Fil_L^tH^i\subset L^pH^i.
\end{equation}
Notice also that
$$
J^{r+t}\subset I^{r+t}\subset J^r+L^t
\qquad
\text{for every $r,t\in\N$}
$$
which implies that the $I$-adic topology on any $A/L^t$-module
agrees with the $J$-adic topology. By the same token, an
$A/L^t$-module is an $(A,\cT_I)$-module of analytically
finite type if and only if it is an $(A,\cT_J)$-module of
analytically finite type. Taking into account proposition
\ref{prop_analyt-noetherian}(i) we deduce that $X$ is also
an analytically noetherian $(A,\cT_J)$-scheme and $\cF/L^t\cF$
is analytically of finite type relative to the topology
$\cT_J$ as well, for every $t\in\N$. Then, theorem
\ref{th_analyt-proper-finiteness} says that
$H^i_t:=H^i(X,\cF/L^t\cF)$ is an $A$-module analytically
of finite type relative to the topology $\cT_J$, hence also
relative to the topology $\cT_I$, since it is an $A/L^tA$-module.
Moreover, by inductive assumption the linear topology on
$H^i_t$ defined by the filtration $\Fil^\bullet_JH^i_t$
agrees with the $J$-adic topology, hence the linear
topology on $H^i_t$ defined by the filtration
$\Fil^\bullet_IH^i_t$ agrees with the $I$-adic topology, by
the foregoing observation.

Now, the image of $\Fil_I^nH^i$ in $H^i_t$ lies in
the intersection of $\Fil_I^nH^i_t$ with the submodule
$H^i/\Fil_L^tH^i\subset H^i_t$, and the $I$-adic topology on
the latter agrees with the one induced by the $I$-adic
topology of $H^i_t$ (corollary \ref{cor_analyt-noether-TARP}
and lemma \ref{lem_Kato-Fuji}(ii)), so
$$
\Img(\Fil_I^nH^i\to H^i/\Fil_L^tH^i)\subset I^p(H^i/\Fil_L^tH^i)
\qquad
\text{for some  $n\in\N$}
$$
{\em i.e.} $\Fil_I^nH^i\subset I^pH^i+\Fil_L^tH^i$, and
combining with \eqref{eq_almost-over}, the assertion
follows.
\end{proof}

\sset\subsubsection{}\label{subsec_coh-of-formal-sch}
Keep the notation of \eqref{subsec_filtration-on-coh},
and define $\pi:X^\wedge\to X$ and $\cF^\wedge$
as in \eqref{subsec_complet-univ-an-noeth-sch}.
Morever, for every $i,n\in\N$ endow $H^i(X,\cF)$ with its
$I$-adic topology and $H^i(X,\cF/I^n\cF)$ with the discrete
topology, and denote by $A^\wedge$ and $H^i(X,\cF)^\wedge$ the
separated completions of $A$ and respectively $H^i(X,\cF)$.

\begin{corollary}\label{cor_rompiballe}
In the situation of \eqref{subsec_coh-of-formal-sch},
the following holds :
\begin{enumerate}
\item
There exist natural $A^\wedge$-linear isomorphisms
$$
H^i(X^\wedge,\cF^\wedge)\isom\lim_{n\in\N}H^i(X,\cF/I^n\cF)
\xleftarrow{\sim}H^i(X,\cF)^\wedge
\qquad
\text{for every $i\in\N$}.
$$
\item
For $i=0$, the maps of {\em (i)} are even isomorphisms
of topological $A^\wedge$-modules.
\item
For every quasi-coherent $\cO_{\!X}$-module $\cG$, there
exists a natural isomorphism
$$
H^i(X^\wedge,\pi^*\cG)\isom A^\wedge\otimes_AH^i(X,\cG)
\qquad
\text{for every $i\in\N$}.
$$
\end{enumerate}
\end{corollary}
\begin{proof}(i): Fix a finite affine covering
$U_\bullet:=(U_\lambda~|~\lambda\in\Lambda)$ of $X$, and
for every $\lambda\in\Lambda$ let $U_\lambda^\wedge$ be
the completion of $U_\lambda$ along its closed subscheme
$\Spec\,A/I\times_{\Spec\,A}U$, . In light of lemma
\ref{lem_left-adj}(i) we have a natural identification
of complexes of topological $A^\wedge$-modules
$$
\bar C{}^\bullet_\mathrm{alt}(U^\wedge_\bullet,\cF^\wedge)\isom
\lim_{n\in\N}\,\bar C{}^\bullet_\mathrm{alt}(U_\bullet,\cF/I^n\cF)
$$
(where $\bar C{}^q_\mathrm{alt}(U_\bullet,\cF/I^n\cF)$ is
endowed with the discrete topology, for every $n,q\in\N$).
On the other hand, by theorems \ref{th_coh-vanish}(i) and
\ref{th_Cech-resolve}(ii), we have natural isomorphisms
$$
H^i(X,\cF/I^n\cF)\isom
H^i\bar C{}^\bullet_\mathrm{alt}(U_\bullet,\cF/I^n\cF)
\qquad
H^i(X^\wedge,\cF^\wedge)\isom
H^i\bar C{}^\bullet_\mathrm{alt}(U^\wedge_\bullet,\cF^\wedge)
$$
for every $i,n\in\N$. Taking into account \cite[Th.3.5.8]{We}
there follows a short exact sequence
\set\begin{equation}\label{eq_Mittag-gut}
0\to\lim_{n\in\N}{}^{\!1}\,H^{i-1}(X,\cF/I^n\cF)\to H^i(X^\wedge,\cF^\wedge)
\to\lim_{n\in\N}\,H^i(X,\cF/I^n\cF)\to 0
\end{equation}
of $A^\wedge$-modules, for every $i\in\N$. Set also
$$
F^i_k:=\Img(H^i(X,\cF)\to H^i(X,\cF/I^k\cF))
\qquad
\text{for every $i,k\in\N$}.
$$
We notice :

\begin{claim}\label{cl_Mittag-gut}
$\displaystyle{\lim_{n\in\N}{}^{\!1}}\,H^{i-1}(X,\cF/I^n\cF)=0$
for every $i\in\N$.
\end{claim}
\begin{pfclaim} By virtue of \cite[Prop.3.5.7]{We}, it suffices
to check that the descending system of submodules
$$
(M^{i-1}_{k,k+n}:=
\Img(H^{i-1}(X,\cF/I^{k+n}\cF)\to H^{i-1}(X,\cF/I^k\cF))~|~n\in\N)
$$
is stationary, for every $i,k\in\N$. However, the short exact
sequence of $\cO_{\!X}$-modules
$$
0\to I^n\cF\to\cF\to\cF/I^n\cF\to 0
$$
induces a commutative diagram of $A$-modules with exact rows
$$
\xymatrix{ H^{i-1}(X,\cF) \ar[r] \ddouble &
H^{i-1}(X,\cF/I^{n+k}\cF) \ar[r] \ar[d] &
H^i(X,I^{n+k}\cF) \ar[r] \ar[d] & H^i(X,\cF) \ddouble \\
H^{i-1}(X,\cF) \ar[r] & H^{i-1}(X,\cF/I^k\cF)
\ar[r] & H^i(X,I^k\cF) \ar[r] & H^i(X,\cF)
}$$
which shows that
$$
F^{i-1}_k\subset M^{i-1}_{k,k+n}
\qquad
\text{for every $i,k,n\in\N$}.
$$
On the other hand, corollary \ref{cor_analyt-finiteness}(ii)
implies that
$$
M^{i-1}_{k,k+n}\subset\Ker(H^{i-1}(X,\cF/I^k\cF)\to H^i(X,I^k\cF))
=F^{i-1}_k
$$
for every $k\in\N$ and every sufficiently large $n\in\N$.
The claim follows.
\end{pfclaim}

From \eqref{eq_Mittag-gut} and claim \ref{cl_Mittag-gut},
we already get the first stated isomorphism. Next, the
proof of claim \ref{cl_Mittag-gut} shows that the system
of inclusion maps $(F^i_n\to H^i(X,\cF/I^n\cF)~|~n\in\N)$
induces an isomorphism of $A^\wedge$-modules
$$
\lim_{n\in\N}\,F^i_n\isom\lim_{n\in\N}\,H^i(X,\cF/I^n\cF)
\qquad
\text{for every $i\in\N$}.
$$
But the system $(F^i_n~|~n\in\N)$ is also isomorphic to
the system $(H^i/\Fil^n_IH^i~|~n\in\N)$ (notation of
\eqref{subsec_filtration-on-coh}), and lastly, from
corollary \ref{cor_analyt-finiteness}(i) we obtain
a natural isomorphism
$$
\lim_{n\in\N}\,H^i/\Fil^n_IH^i\isom H^i(X,\cF)^\wedge
$$
whence the second stated isomorphism.

(ii): The topology of $H^0(X,\cF)^\wedge$ is $I$-adic,
by remark \ref{rem_completion-of-topring}(ii,iv).
On the other hand, the topology of $H^0(X^\wedge,\cF^\wedge)$
is the linear topology given by the system of submodules
$$
(H^0(X^\wedge,(I^n\cF)^\wedge)~|~n\in\N)
$$
by remark \ref{rem_top-on-sections}. Denote $H^0(X,I^n\cF)^\wedge$
the $I$-adic completion of $H^0(X,I^n\cF)$, for every $n\in\N$;
by (i) we get a commutative diagram of $A^\wedge$-modules
\set\begin{equation}\label{eq_compare-tops}
{\diagram
H^0(X^\wedge,(I^n\cF)^\wedge) \ar[r] \ar[d] &
H^0(X,I^n\cF)^\wedge \ar[d] \\
H^0(X^\wedge,\cF^\wedge) \ar[r] & H^0(X,\cF)^\wedge
\enddiagram}
\end{equation}
whose horizontal arrows are isomorphisms, and whose left
(resp. right) vertical arrow is induced by the inclusion
$(I^n\cF)^\wedge\subset\cF^\wedge$ (resp. is the $I$-adic
completion of the natural map $H^0(X,I^n\cF)\to H^0(X,\cF)$).
By corollary \ref{cor_analyt-finiteness}(i), the topology
on $H^0(X,\cF)^\wedge$ given by the filtration
$(H^0(X,I^n\cF)^\wedge~|~n\in\N)$ agrees with the topology
given by the filtration $((I^nH^0(X,\cF))^\wedge~|~n\in\N)$,
and the latter is none else than the $I$-adic filtration,
again by remark \ref{rem_completion-of-topring}(ii,iv).
We conclude already that the bottom horizontal arrow of
\eqref{eq_compare-tops} is an isomorphism of topological
$A^\wedge$-modules. The same holds for the natural map
$H^0(X^\wedge,\cF^\wedge)\to\lim_{n\in\N}H^0(X,\cF/I^n\cF)$,
by the discussion of \eqref{subsec_complet-along}.

(iii): Since $X$ is quasi-compact and quasi-separated,
$\cG$ is the colimit of the filtered family
$(\cG_\lambda~|~\lambda\in\Lambda)$ of its quasi-coherent
$\cO_{\!X}$-submodules of finite type (proposition
\ref{prop_fp-approx}), hence $\pi^*\cG$ is the colimit
of the system $(\pi^*\cG_\lambda~|~\lambda\in\Lambda)$.
Since $X^\wedge$ is a spectral topological space, proposition
\ref{prop_dir-im-and-colim} then reduces to showing the
assertion for each $\cG_\lambda$, so we may assume from start
that $\cG$ is analytically of finite type, in which case it
suffices to apply (i), propositions
\ref{prop_like-noetherian}(i) and \ref{prop_up-to-completion}(i),
and theorem \ref{th_analyt-proper-finiteness}.
\end{proof}

\begin{remark}\label{rem_conclusion}
(i)\ \
The assertions of theorem \ref{th_analyt-proper-finiteness}
and of its corollaries \ref{cor_analyt-finiteness} and
\ref{cor_rompiballe} hold also in case $X$ is projective
but not necessarily finitely presented : indeed, if
$i:X\to\P^n_A$ is a closed immersion in a projective
space over $\Spec\,A$, then
$H^\bullet(X,\cF)=H^\bullet(\P^n_A,i_*\cF)$, and $i_*\cF$
is clearly an $\cO_{\P^n_A}$-module of analytically finite
type, so we are reduced to the case where $X=\P^n_A$, which
is covered by theorem \ref{th_analyt-proper-finiteness}.

(ii)\ \
More generally, theorem \ref{th_analyt-proper-finiteness}
and its two  corollaries hold for any proper $A$-scheme $X$
(and any $\cO_{\!X}$-module $\cF$ of analytically finite
type) : indeed, according to \cite[Th.4.1 and Th.4.3]{Con2}
there exists a locally closed immersion $i:X\to X'$ of
$A$-scheme, with $X'$ proper and finitely presented over
$\Spec\,A$. Since $X$ is proper, $i$ is even a closed
immersion (\cite[Ch.II, Cor.5.4.3]{EGAII}), so again it
suffices to apply the theorem to the $\cO_{\!X'}$-module
$i_*\cF$.

(iii)\ \
Theorem \ref{th_analyt-proper-finiteness} and corollary
\ref{cor_rompiballe} have also been announced, respectively
as theorem C.3.1 and theorem C.3.3, in K.Fujiwara and F.Kato's
treatise \cite{Fu-Ka}; however, their proofs have been
postponed to the second part of this work, which has not
appeared yet.
\end{remark}

\begin{corollary} Let $A$ be a universally analytically
noetherian ring, $\bff:=(f_1,\dots,f_r)$ a sequence of
elements of $A$ that generates an ideal of adic definition.
Then $A$ satisfies condition $\mathrm{(a)^{un}_\bff}$ of
\eqref{subsec_badabum}.
\end{corollary}
\begin{proof} Let $X:=\Spec\,A$, denote by $I\subset A$
the ideal generated by $\bff$, and by $\cI\subset\cO_{\!X}$
the quasi-coherent ideal arising from $I$. Let also
$\pi:Y\to X$ be the blowing up morphism of $\cI$. It
follows easily from remark \ref{rem_blowing-up}(iii)
that the support of the $A$-module $H^p:=H^p(Y,\cO_Y)$
is contained in $\Spec\,A/I$, for every $p>0$, and the
same holds for the kernel and cokernel of the natural
map $\phi:A\to H^0(Y,\cO_Y)$. Combining with theorem
\ref{th_analyt-proper-finiteness} and remark
\ref{rem_anal-noetherian}(iii), we deduce that there
exists $n\in\N$ such that $I^n\cdot H^p=0$ for every $p>0$,
and $I^n\cdot\Ker\,\phi=I^n\cdot\Coker\,\phi=0$. Then it
suffices to invoke theorem \ref{th_too-difficult-for-me}.
\end{proof}

\subsection{Continuous valuations}
\label{sec_affinoid-rings}
As already mentioned in \eqref{subsec_not-pro-constr},
the study of the valuation spectrum of a topological
ring $A$ will lead us to consider certain subsets of
$\Spv\,A$ that are not pro-constructible, but nevertheless
are spectral spaces, with the topology induced via the
inclusion map into $\Spv\,A$.

\sset\subsubsection{}\label{subsec_I-admissible}
As a first step, let $A$ be any ring, $I\subset A$
a finitely generated ideal; we attach to $I$ a set of
specializations $S_I$ of $\Spv\,A$ (see
\eqref{subsec_not-pro-constr}), by declaring that
$(v,w)\in S_I$ if and only if :
\begin{itemize}
\item
either, $w(I)\neq\{0\}$ and $w$ is a primary specialization
of $v$
\item
or else, $v=w$.
\end{itemize}
The $S_I$-admissible specializations in $\Spv\,A$ will
be simply called {\em $I$-admissible}, and the $S_I$-closed
subsets of $\Spv\,A$ shall be called {\em $I$-closed}. It
is clear that the set of $I$-admissible specializations of
$\Spv\,A$ is transitive, and it satisfies condition (S3) of
corollary \ref{cor_Huber-crit}, by lemma
\ref{lem_Spezializerungen}(i). We shall show hereafter
that it fulfills also conditions (S1) and (S2) of
proposition \ref{prop_Huber-crit}. This shall be achieved
in several steps; to start out, for any abelian ordered
group $\Delta$, let us say that an element $\delta\in\Delta_\circ$
is {\em final in $\Delta$} if the following holds. For every
$\gamma\in\Delta$ there exists $n\in\N$ such that
$\delta^n<\gamma$. We remark :

\begin{lemma}\label{lem_I-admit}
In the situation of \eqref{subsec_I-admissible},
let $v$ be a valuation of $A$ with
$$
v(I)\neq\{0\}
\qquad\text{and}\qquad
v(I)\cap c\Gamma_{\!v}=\emptyset.
$$
Then we have :
\begin{enumerate}
\item
There exists a smallest convex subgroup $\Delta$ of\/
$\Gamma_{\!v}$ with $v(I)\cap\Delta\neq\emptyset$.
\item
$c\Gamma_{\!v}\subset\Delta$.
\item
For every $i\in I$, either $v(i)<\delta$ for every
$\delta\in\Delta$, or else $v(i)$ is final in $\Delta$.
\end{enumerate}
\end{lemma}
\begin{proof}(i): Let $a_1,\dots,a_n$ be a finite system
of generators of $I$; the stated assumptions on $v$
imply that $v(a_i)<1$ for every $i=1,\dots,n$, and
$v(a_k)>0$ for at least one index $k\leq n$. Let $\Delta$
be the convex hull of the subgroup $\Delta_0$ of $\Gamma_{\!v}$
generated by $\max\{v(a_i)~|~i=1,\dots,n\}$. It remains
only to check that $\Delta\subset\Delta'$ for every
convex subgroup $\Delta'\subset\Gamma_{\!v}$ with
$\Delta'\cap v(I)\neq\emptyset$. To this aim, notice
first that $c\Gamma_{\!v}\subset\Delta'$ for every such
$\Delta'$ : indeed, otherwise there exists
$\gamma\in c\Gamma_{\!v}$ such that $\gamma<\delta$ for
every $\delta\in\Delta'$. Then, pick $i\in I$ with
$v(i)\in\Delta'$; it follows that $1>v(i)>\gamma$, so
$v(i)\in c\Gamma_{\!v}$, contradicting our assumptions.
Now, let $i\in I$ be any element with $v(i)\in\Delta'$,
and write $i=a_1b_1+\cdots+a_nb_n$ for some
$b_1,\dots,b_n\in A$; since $c\Gamma_{\!v}\subset\Delta'$
and $v(i)\leq\max\{v(a_i)\cdot v(b_i)~|~i=1,\dots,n\}$, it
follows that $v(a_k)\in T$ for some $k\leq n$ (details
left to the reader). Then clearly $\Delta_0\subset\Delta'$,
and the contention follows. Moreover, the argument also
shows (ii).

(iii): Let $a_0\in I$ be any element; the sequence
$a_0,a_1,\dots,a_n$ is also a system of generators
for $I$, and the foregoing shows that $\Delta$ is
the convex hull of the subgroup generated by
$\max\{v(a_i)~|~i=0,\dots,n\}$. The assertion follows
immediately.
\end{proof}

\begin{definition}\label{def_cGamma-I}
Keep the situation of \eqref{subsec_I-admissible}.
To every $v\in\Spv\,A$ we attach a convex subgroup
$$
c\Gamma_{\!v}(I)\subset\Gamma_{\!v}
$$
as follows :
\begin{itemize}
\item
If $v(I)=\{0\}$, we set $c\Gamma_{\!v}(I):=\Gamma_{\!v}$.
\item
If $v(I)\cap c\Gamma_{\!v}\neq\emptyset$, we set
$c\Gamma_{\!v}(I):=c\Gamma_{\!v}$.
\item
Otherwise, $v$ fulfills the conditions of lemma
\ref{lem_I-admit}, and we let $c\Gamma_{\!v}(I)$ be
the smallest of the convex subgroups $\Delta$ such that
$v(I)\cap\Delta\neq\emptyset$.
\end{itemize}
\end{definition}

\begin{lemma}\label{lem_cGamma-I}
Keep the notation of definition {\em\ref{def_cGamma-I}},
and let $a_1,\dots,a_n$ be a system of generators for the
ideal $I$. We have :
\begin{enumerate}
\item
The following conditions are equivalent :
\begin{enumerate}
\item
$\Gamma_{\!v}=c\Gamma_{\!v}(I)$.
\item
$\Gamma_{\!v}=c\Gamma_{\!v}$, or else $v(a)$ is final in
$\Gamma_{\!v}$ for every $a\in I$.
\item
$\Gamma_{\!v}=c\Gamma_{\!v}$, or else $v(a_i)$ is final in
$\Gamma_{\!v}$, for every $i=1,\dots,n$.
\end{enumerate}
\item
An element $w\in\Spv\,A$ is an $I$-admissible specialization
of $v$ if and only if  $w=v^\Delta$ for a convex subgroup
$\Delta$ of\/ $\Gamma_{\!v}$ with $c\Gamma_{\!v}(I)\subset\Delta$.
\item
$v$ has no proper $I$-admissible specializations if and only
if\/ $\Gamma_{\!v}=c\Gamma_{\!v}(I)$.
\end{enumerate}
\end{lemma}
\begin{proof}(i): The equivalence of (a) and (b) follows
directly from the definitions and lemma
\ref{lem_I-admit}(iii). The equivalence of (b) and (c)
follows by remarking that the set
$$
\{a\in A~|~
\text{$v(a)=0$ or else $v(a)$ is final in $\Gamma_{\!v}$}\}
$$
is an ideal in $A$, provided $\Gamma_{\!v}\neq c\Gamma_{\!v}$
(details left to the reader).

(ii) and (iii) are immediate from the definitions.
\end{proof}

\begin{lemma}\label{lem_characterize-I-closed}
Let $A$ be a ring, $I\subset A$ a finitely generated
ideal. We have :
\begin{enumerate}
\item
Let also $f_\bullet:=(f_0,f_1,\dots,f_n)$ be a finite
sequence of elements of $A$, and $J\subset A$ the ideal
generated by $f_\bullet$. Suppose that the radical of $J$
contains $I$. Then the rational subset
$R_A\bigl(\frac{f_1}{f_0},\dots,\frac{f_n}{f_0}\bigr)$ is
$I$-closed.
\item
Let $v$ be any valuation of $A$ such that
$\Gamma_{\!v}=c\Gamma_{\!v}(I)$. Then for every open neighborhood
$U$ of $v$ in $\Spv\,A$ there exists a finite sequence
$f_\bullet:=(f_0,f_1,\dots,f_n)$ of elements of $A$ such that :
\begin{enumerate}
\item
$v\in R_A\bigl(\frac{f_1}{f_0},\dots,\frac{f_n}{f_0}\bigr)\subset U$.
\item
$I$ is contained in the radical of the ideal generated
by $f_\bullet$.
\end{enumerate}
\end{enumerate} 
\end{lemma}
\begin{proof}(i): Let $w\in\Spv\,A$ be an $I$-admissible
specialization of some
$v\in R:=R_A\bigl(\frac{f_1}{f_0},\dots,\frac{f_n}{f_0}\bigr)$.
If $v(I)=\{0\}$, then $v=w$, and there is nothing to
prove. Thus, we may assume $v(I)\neq\{0\}$, in which
case $w(I)\neq\{0\}$ as well. Since $w$ is a primary
specialization of $v$, we have $w(f_i)\leq w(f_0)$ for
every $i=1,\dots,n$. Suppose, by way of contradiction,
that $w\notin R$; then we must have $w(f_0)=0$, and
consequently $w(f_i)=0$ for every $i\leq n$ as well.
Thus, $J\subset\Ker\,w$, and since the support of $w$
is a prime ideal, it contains also the radical of $J$;
hence, $I\subset\Ker\,w$, a contradiction.

(ii): We may assume that
$U=R_A\bigl(\frac{g_1}{g_0},\dots,\frac{g_k}{g_0}\bigr)$ for
some $g_0,g_1,\dots,g_k\in A$. Pick also a finite system
$a_1,\dots,a_n$ of generators of $I$. We distinguish two
cases :

$\bullet$ Suppose first that $\Gamma_{\!v}=c\Gamma_{\!v}$.
In this case, since $v(g_0)\neq 0$, there exists $d\in A$
such that $v(g_0)\cdot v(d)=v(g_0d)>1$. Then
$v\in R_A\bigl(\frac{1}{g_0d}\bigr)\cap U$, and we may take
$f_\bullet:=(g_0d,\dots,g_kd,1)$.

$\bullet$ Lastly, suppose that $\Gamma_{\!v}\neq c\Gamma_{\!v}$.
In this case, since $v(g_0)\neq 0$, lemma \ref{lem_cGamma-I}(i)
implies that there exists $r\in\N$ with
$v(g_0)\geq v(a_i^r)$ for every $i=1,\dots,n$. Then
$v\in R_A\bigl(\frac{a_1^r}{g_0},\dots,\frac{a_n^r}{g_0}\bigr)\cap U$,
and the sequence $f_\bullet:=(g_0,\dots,g_k,a_1^r,\dots,a_n^r)$
will do.
\end{proof}

\begin{definition}\label{def_Spv-I}
Let $A$ be any ring and $I\subset A$ any finitely
generated ideal.
\begin{enumerate}
\item
We set
$$
\Spv(A,I):=\{v\in\Spv\,A~|~\Gamma_{\!v}=c\Gamma_{\!v}(I)\}
$$
and we endow $\Spv(A,I)$ with the topology induced by
$\Spv\,A$ via the inclusion map.
\item
The {\em rational subsets} of $\Spv(A,I)$ are the subsets
of the form
$$
\Spv(A,I)\cap R_A\Bigl(\frac{f_1}{f_0},\dots,\frac{f_n}{f_0}\Bigr)
$$
where $f_\bullet:=(f_0,f_1,\dots,f_n)$ is any sequence of
elements of $A$ such that $I$ is contained in the radical
of the ideal generated by $f_\bullet$.
\end{enumerate}
\end{definition}

\begin{theorem}\label{th_Spv-I}
With the notation of definition {\em\ref{def_Spv-I}},
we have :
\begin{enumerate}
\item
The topological space $\Spv(A,I)$ is spectral, and the
retraction
$$
r_I:\Spv\,A\to\Spv(A,I)
\qquad
v\mapsto v^{c\Gamma_{\!v}(I)}
$$
is spectral.
\item
The rational subsets of\/ $\Spv(A,I)$ are constructible
in $\Spv(A,I)$, and form a basis of the topology of\/
$\Spv(A,I)$ that is closed under finite intersections.
\item
A subset $T\subset\Spv(A,I)$ is constructible in $\Spv(A,I)$
if and only if\/ $r^{-1}_IT$ is constructible in $\Spv\,A$.
\end{enumerate}
\end{theorem}
\begin{proof}(i): From lemma \ref{lem_characterize-I-closed}
we see that the set of $I$-admissible specializations of
$\Spv\,A$ fulfills conditions (S1) and (S2) of proposition
\ref{prop_Huber-crit}, whence the assertion.

(ii): Lemma \ref{lem_characterize-I-closed} already
shows that the rational subsets form a basis $\cB$ of
the topology of $\Spv(A,I)$. If
$R_A\bigl(\frac{f_1}{f_0},\dots,\frac{f_n}{f_0}\bigr)$ and
$R_A\bigl(\frac{g_1}{g_0},\dots,\frac{g_m}{g_0}\bigr)$ are
any two rational open subset of $\Spv\,A$ such that
the radicals of the ideals $J$ and $J'$ generated
respectively by $f_\bullet$ and $g_\bullet$ contain $I$,
then the same holds for the radical of the ideal $JJ'$;
taking into account remark \ref{rem_Spv-of-ring}(i), it
follows that $\cB$ is closed under finite intersections.

Lastly, let $f_\bullet$ and $J$ be as in the foregoing,
with $I$ contained in the radical of $J$, and set
$R:=R_A\bigl(\frac{f_1}{f_0},\dots,\frac{f_n}{f_0}\bigr)$.
By lemma \ref{lem_characterize-I-closed}, the set
$R$ is $I$-closed, so $r_I(R)=R\cap\Spv(A,I)$; since
$R$ is quasi-compact in $\Spv\,A$, the subset $r_I(R)$
is quasi-compact in $\Spv(A,I)$, so $R\cap\Spv(A,I)$
is a constructible open subset of $\Spv(A,I)$.

(iii) follows from proposition \ref{prop_Huber-crit}(iii).
\end{proof}

\begin{example}\label{ex_Spv-A-I}
Let $A$ be any ring.

(i)\ \
We may take $I=0$, in which case $c\Gamma_{\!v}(I)=\Gamma_{\!v}$
for every $v\in\Spv\,A$, so that $\Spv(A,0)=\Spv\,A$. Also,
the rational subsets of $\Spv(A,0)$ are the same as the
rational subsets of $\Spv\,A$.

(ii)\ \
We may also take $I=A$, in which case
$c\Gamma_{\!v}(I)=c\Gamma_{\!v}$ for every $v\in\Spv\,A$, so
that $\Spv(A,A)=\{v\in\Spv\,A~|~\Gamma_{\!v}=c\Gamma_{\!v}\}$.
The rational subsets of $\Spv(A,A)$ can be described as
the sets of the form
$$
\Spv(A,A)\cap\{v\in\Spv\,A~|~
               v(f_1)\leq v(f_0),\dots,v(f_n)\leq v(f_0)\}
$$
where $f_0,f_1,\dots,f_n$ is any sequence of elements
with $\sum_{i=1}^nf_iA=A$ (details left to the reader).
\end{example}

\begin{lemma}\label{lem_was-rem-Spv-A-I}
Let $A$ be any ring, and $I,J\subset A$ two finitely
generated ideals. We have :
\begin{enumerate}
\item
If $I$ is contained in the radical of $J$, then
$\Spv(A,J)\subset\Spv(A,I)$. The inclusion map
$\Spv(A,J)\to\Spv(A,I)$ is not spectral in general.
On the other hand, there is a spectral retraction
$$
r_{I,J}:\Spv(A,I)\to\Spv(A,J)
\qquad\text{such that}\qquad
r_{I,J}\circ r_I=r_J.
$$
\item
Especially, the spectral map $r_{I,I+J}$ is well defined,
and it coincides with the restriction of $r_J$ to $\Spv(A,I)$.
\item
$\Spv(A,I\cdot J)=\Spv(A,I)\cup\Spv(A,J)$.
\item
$\Spv(A,I+J)=\Spv(A,I)\cap\Spv(A,J)$.
\end{enumerate}
\end{lemma}
\begin{proof}(i): Indeed, the assumption easily implies
that $c\Gamma_{\!v}(J)\subset c\Gamma_{\!v}(I)$ for every
$v\in\Spv\,A$, from which the assertion follows immediately.

(ii): Let $v\in\Spv\,A$ be any element, and set
$w:=r_I(v)$, $u:=r_J(w)$; then $u$ does not have proper
$J$-admissible specializations, nor any proper $I$-admissible
specialization, because any $I$-admissible specialization
of $u$ would be also an $I$-admissible specialization of
$w$. It follows easily that $u$ is a primary specialization
of $v$ that does not have any proper $(I+J)$-specialization,
whence the assertion. 

\begin{claim}\label{cl_about-plus-and-dot}
For every $v\in\Spv\,A$ we have :
$$
c\Gamma_{\!v}(I+J)=c\Gamma_{\!v}(I)\cap c\Gamma_{\!v}(J)
\qquad\text{and}\qquad
c\Gamma_{\!v}(I\cdot J)=c\Gamma_{\!v}(I)\cup c\Gamma_{\!v}(J).
$$
\end{claim}
\begin{pfclaim} We consider first the assertion for
$c\Gamma_{\!v}(I\cdot J)$ : it is easily seen that
$v(IJ)=\{0\}$ if and only if either $v(I)=\{0\}$ or
$v(J)=\{0\}$, so the assertion holds in case $v(IJ)=\{0\}$.
Similarly, $v(IJ)\cap c\Gamma_{\!v}\neq\emptyset$ if and
only if both $v(I)$ and $v(J)$ intersect $c\Gamma_{\!v}$,
so the assertion holds also in case $v(IJ)$ intersects
$c\Gamma_{\!v}$. Lastly, suppose that $v(IJ)\neq\{0\}$
and $v(IJ)\cap c\Gamma_{\!v}=\emptyset$, and pick
a finite system of generators $a_1,\dots,a_n$ (resp.
$b_1,\dots,b_m$) for $I$ (resp. for $J$), with
$$
v(a_1)=\max(v(a_i)~|~i=1,\dots,n)
\qquad\text{and}\qquad
v(b_1)=\max(v(b_i)~|~i=1,\dots,m).
$$
Then the system $(a_ib_j~|~i\leq n,\ j\leq m)$ generates
$IJ$, and $v(ab)=\max(v(a_ib_j)~|~i\leq n,\ j\leq m)$.
Hence $v(a_1)$ (resp. $v(b_1)$, resp. $v(a_1b_1)$) is
final in $c\Gamma_{\!v}(I)$ (resp. in $c\Gamma_{\!v}(J)$,
resp. in $c\Gamma_{\!v}(IJ)$) by lemma \ref{lem_I-admit}(iii).
Without loss of generality we may assume that
$c\Gamma_{\!v}(J)\subset c\Gamma_{\!w}(I)$, in which case
there exists $n\in\N$ such that $v(a_1^n)\leq v(b_1)$,
hence $v(a_1b_1)^n\geq v(a_1)^{2N}$, which shows that
$c\Gamma_{\!v}(IJ)\subset c\Gamma_{\!v}(I)$, and the
converse inclusion is clear.

We argue similarly for the calculation of
$c\Gamma_{\!v}(I+J)$ : first, it is clear that $v(I+J)=\{0\}$
if and only if $v(I)=v(J)=\{0\}$. Next, we have
$v(I+J)\cap c\Gamma_{\!v}\neq\emptyset$ if and only
if $(v(I)\cup v(J))\cap c\Gamma_{\!v}\neq\emptyset$
(details left to the reader). Laslty, if $v(IJ)\neq\{0\}$
and $v(IJ)\cap c\Gamma_{\!v}=\emptyset$, pick $a_1,\dots,a_n$
and $b_1,\dots,b_m$ as in the foregoing; without loss
of generality we may assume that $v(a_1)\geq v(b_1)$.
Since the system $a_1,\dots,a_n,b_1,\dots,b_m$ generates
$I+J$, we deduce that $v(a_1)$ is final in $c\Gamma_{\!v}(I+J)$
(lemma \ref{lem_I-admit}(iii)), whence the contention.
\end{pfclaim}

Assertion (iii) and (iv) follow immediately from claim
\ref{cl_about-plus-and-dot} : details left to the reader.
\end{proof}

\begin{example} Let $A$ be any ring, $I\subset A$ a finitely
generated ideal.

(i)\ \
Every element of the subset $L:=\{v\in\Spv\,A~|~v(I)=\{0\}\}$
has no proper $I$-admissible specializations, so $L$ lies
in $\Spv(A,I)$. Clearly, $L$ is constructible in $\Spv\,A$,
and $r^{-1}_I(L)=L$, so $L$ is constructible in $\Spv(A,I)$,
by theorem \ref{th_Spv-I}(iii).

(ii)\ \
Let $(\Spv\,A)_0$ be the set of trivial valuations of $A$;
from remark \ref{rem_Spv-of-ring}(iii) we know already that
$(\Spv\,A)_0$ is a spectral space, with the topology induced
by the inclusion into $\Spv\,A$. Moreover, if
$R:=R_A\bigl(\frac{f_1}{f_0},\dots,\frac{f_n}{f_0}\bigr)$ is
any rational subset of $\Spv(A,I)$, then
$R\cap(\Spv\,A)_0=\{v\in(\Spv\,A)_0~|~v(f_0)=1\}$, which
is a constructible (open) subset of $(\Spv\,A)_0$. Hence,
the inclusion map $(\Spv\,A)_0\to\Spv(A,I)$ is spectral.

(iii)\ \
Pick a finite system $a_1,\dots,a_n$ of generators for $I$,
and set
$$
\begin{aligned}
E_0:=\, & \{v\in(\Spv\,A)_0~|~v(I)=\{0\}\} \\
E_i:=\, & \{v\in\Spv\,A~|~
\text{$v(a_i)=1$ and $v(a)\leq 1$ for every $a\in A$}\}
& \text{for every $i=1,\dots,n$}.
\end{aligned}
$$
Then $E_0,\dots,E_n$ are pro-constructible subsets of
$\Spv\,A$ that are closed under $I$-admissible
specializations and generizations. We have
$(\Spv\,A)_0=\Spv(A,I)\cap\bigcup_{i=0}^nE_i$, and
therefore $r^{-1}_I(\Spv\,A)_0=\bigcup_{i=0}^nE_i$,
by corollary \ref{cor_Huber-crit}.
\end{example}

We are now ready to introduce the {\em spectrum of
continuous valuations} of a topological ring.

\begin{definition}\label{def_continuous-vals}
Let $(A,\cT_A)$ be any topological ring, and
$v$ a valuation of $A$.
\begin{enumerate}
\item
We endow $\Gamma_{\!v\circ}$ with a topology $\cT_{\Gamma_{\!v}}$, by
ruling that a subset $U\subset\Gamma_{\!v\circ}$ is open if
either $0\notin U$, or else there exists $\delta\in\Gamma_{\!v}$
such that $\{\gamma\in\Gamma_{\!v\circ}~|~\gamma<\delta\}\subset U$.
Then, we say that $v$ is {\em continuous} if it is a
continuous map $(A,\cT_A)\to(\Gamma_{\!v\circ},\cT_{\Gamma_{\!v}})$.
\item
We set
$$
\Cont(A):=\{v\in\Spv\,A~|~\text{$v$ is continuous}\}
\qquad
\Cont^+(A):=\Cont(A)\cap\Spv^+(A)
$$
and we endow $\Cont(A)$ and $\Cont^+(A)$ with the topology
induced by $\Spv\,A$.
\end{enumerate}
\end{definition}

\begin{remark}\label{rem_Cont-A}
Let $A$ be any topological ring, any $v$ any valuation
of $A$.

(i)\ \
Endow $\kappa(v)$ with its valuation topology $\cT_v$
(see definition \ref{def_valuation-toplog}(i)).
A simple inspection of the definitions shows that the
residual valuation $\bar v:(\kappa(v),\cT_v)\to\Gamma_{\!v\circ}$
of $v$ is continuous. Moreover, $v$ is continuous if
and only if the same holds for the natural map
$\pi_v:A\to(\kappa(v),\cT_v)$. Indeed, clearly if $\pi_v$
is continuous, the same holds for $v=\bar v\circ\pi_v$.
Conversely, suppose that $v$ is continuous; it suffices
to check that $\pi_v$ is continuous at the point $0\in A$.
However, for every $\delta\in\Gamma_{\!v}$ let
$U_\delta:=\{\gamma\in\Gamma_{\!v\circ}~|~\gamma<\delta\}$;
then the family $(U_\delta~|~\delta\in\Gamma_{\!v})$ is a
fundamental system of open neighborhoods of
$0\in\Gamma_{\!v\circ}$, and
$(\bar v{}^{-1}U_\delta~|~\delta\in\Gamma_{\!v})$ is a
fundamental system of open neighborhoods of $0\in\kappa(v)$,
whence the contention.

(ii)\ \
It follows easily from (i) that the valuation $v$
is continuous if and only if it is continuous at the
point $0\in A$.

(iii)\ \
Let $\phi:\Gamma\to\Gamma'$ be a morphism of ordered groups,
and endow $\Gamma_{\!\circ}$ and  $\Gamma'_{\!\circ}$ with the
topologies $\cT_\Gamma$ and $\cT_{\Gamma'}$ described in definition
\ref{def_continuous-vals}(i). In general, the induced map
$\phi_\circ:\Gamma_{\!\circ}\to\Gamma'_{\!\circ}$ is {\em not}
necessarily continuous for these topologies; for instance :
\begin{itemize}
\item
If $\Gamma$ is a subgroup of $\Gamma'$, and $\phi$ is the
inclusion map, then $\phi_\circ$ is continuous if and only
if the convex hull of $\Gamma$ equals $\Gamma'$ (see remark
\ref{rem_ordered-gps}(iii)), and if the latter condition
holds, then the topology $\cT_\Gamma$ agrees with the topology
induced by $\cT_{\Gamma'}$.
\item
If $\Gamma'=\Gamma/\Delta$ for a proper convex subgroup
$\Delta\subset\Gamma$, and $\phi$ is the projection, then
$\phi_\circ$ is continuous.
\item
If $\Gamma\neq\{1\}$ and $\Gamma'=\{1\}$, then $\phi_\circ$
is not continuous.
\end{itemize}
Especially, a secondary generization $w$ of $v$ is continuous
if either $v=w$ or else $\Gamma_{\!w}\neq\{1\}$, but $w$ may
fail to be continuous in case $\Gamma_{\!v}\neq\{1\}$ and
$\Gamma_{\!w}=\{1\}$. On the other hand, if $\phi_\circ\circ v$
is a continuous map, then it is easily seen that $v$ is
continuous as well.

(iv)\ \
Let $\Delta\subset\Gamma_{\!v}$ be any convex subgroup
containing $c\Gamma_{\!v}$, denote by
$\Gamma^+_{\!v\circ}\cdot\Delta$ the smallest submonoid
of $\Gamma_{\!v\circ}$ containing $\Gamma^+_{\!v\circ}$ and
$\Delta$, and endow $\Gamma^+_{\!v\circ}\cdot\Delta$ with
the topology $\cT_{\Gamma^+_{\!v\circ}\cdot\Delta}$ induced by
$\cT_{\Gamma_{\!v}}$ via the inclusion map. Let also $\cT_\Delta$
be the topology on $\Delta_\circ$ described in definition
\ref{def_continuous-vals}(i). Then the inclusion map
$\Delta_\circ\to\Gamma^+_{\!v\circ}\cdot\Delta$ admits a
unique continuous retraction
$$
\rho:(\Gamma^+_{\!v\circ}\cdot\Delta,\cT_{\Gamma^+_{\!v\circ}\cdot\Delta})
\to(\Delta_\circ,\cT_\Delta)
\qquad
\text{such that $\rho(\gamma)=0$ whenever
$\gamma\not\in\Delta_\circ$}. 
$$
On the other hand, since $c\Gamma_{\!v}\subset\Delta$, the
map $v$ factors through a (unique) mapping
$v':A\to\Gamma^+_{\!v\circ}\cdot\Delta$, and if $v$ is
continuous, the same holds for $v'$. Clearly
$\rho\circ v'=v^\Delta$, so the latter is continuous
as well. Combining with (iii) and proposition
\ref{prop_reverse-order}, we conclude that every
specialization of a continuous valuation is continuous.
However, in general $\Cont(A)$ is not a closed subset
(not even a pro-constructible subset) of $\Spv\,A$.

(v)\ \
Since $\Spv\,A$ is a $T_0$ topological space, the same
holds for $\Cont(A)$. Especially, the specializations
in $\Cont(A)$ define a partial ordering on $\Cont(A)$ :
see remark \ref{rem_specialize}(ii).

(vi)\ \
It follows also from (ii) that if $\Ker\,v$ is an
open prime ideal, then $v$ is continuous. We shall
say that $v$ is {\em non-analytic} if $\Ker\,v$ is
an open prime ideal. If $v$ is continuous and its
support is not an open prime ideal, we shall say
that $v$ is {\em analytic}. We denote by
$$
\Cont(A)_\mathrm{a}
\qquad\text{and}\qquad
\Cont(A)_\mathrm{na}
$$
the sets of all analytic and respectively non-analytic
valuations of $A$.

(vii)\ \
Let $f:A\to B$ be any continuous ring homomorphism of
topological rings. Then $\Spv(f)$ restricts to continuous
maps
$$
\Cont(f):\Cont(B)\to\Cont(A)
\qquad
\Cont^+(f):\Cont^+(B)\to\Cont^+(A).
$$
Indeed, let $w$ be a continuous valuation of $B$; obviously
$v:=w\circ f:A\to\Gamma_{\!w\circ}$ is a continuous map.
The value group $\Gamma_{\!v}$ of $v$ is a subgroup of
$\Gamma_{\!w}$, and $v$ factors uniquely through a valuation
$v':A\to\Gamma_{\!v}$ that is obviously equivalent to $v$;
by (iii), the continuity of $w$ implies the continuity of
$v'$, whence the claim. It is then also clear that
$\Cont(f)$ in turns restricts to a mapping
$$
\Cont(f)_\mathrm{na}:\Cont(B)_\mathrm{na}\to\Cont(A)_\mathrm{na}.
$$
\end{remark}

\begin{lemma}\label{lem_Cont-A}
Let $A$ be any f-adic topological ring, and $v\in\Cont(A)$. We have :
\begin{enumerate}
\item
$\Cont(A)=\{w\in\Spv\,A~|~
\text{$w(a)$ is final in $\Gamma_{\!w}$, for every
$a\in A^{\circ\circ}$}\}$.
\item
Suppose that $v$ is analytic, and endow $\kappa(v)$ with its
valuation topology $\cT_v$. Then $v$ is a Tate valuation (see
definition {\em\ref{def_valuation-toplog}(iv)}), and the natural
map
$$
\pi_v:A\to(\kappa(v),\cT_v)
$$
is an f-adic continuous ring homomorphism.
\item
Let $w\in\Cont(A)$ be any proper primary specialization of $v$.
Then $w$ is non-analytic.
\item
Let $w\in\Cont(A)_\mathrm{a}$ be a specialization of $v$. Then
$v\in\Cont(A)_\mathrm{a}$ and $w$ is a secondary specialization of $v$.
\item
If $v$ is analytic, the totally ordered set of continuous
generizations of\/ $v$ admits a maximal element, which
is a rank one analytic valuation.
\end{enumerate}
\end{lemma}
\begin{proof} Let $A_0\subset A$ be a subring of definition,
$I\subset A_0$ a finitely generated ideal of adic definition,
and $a_1,\dots,a_n$ a finite system of generators of $I$; for
every $w\in\Spv\,A$ we set $\gamma_w:=\max(w(a_1),\dots,w(a_n))$

(i): The condition is clearly necessary. Conversely, suppose
that $w\in\Spv\,A$ fulfills this condition; by remark
\ref{rem_Cont-A}(ii), it suffices to check that, for every
$\delta\in\Gamma_{\!w}$ there exists an open neighborhood
$U_\delta$ of $0$ in $A$ such that $w(a)<\delta$ for every
$a\in U_\delta$. However, we have $w(a)<1$ for every $a\in I$,
since $I\subset A^{\circ\circ}$, and we may find $k\in\N$
such that $\delta>\gamma^k_w$. Hence $U_\delta:=I^{k+1}$ will do.

(ii): Denote by $\bar v$ the residual valuation of $v$; if
$v$ is analytic, the value $\gamma_v$ lies in $\Gamma_{\!v}$,
and if $x\in\kappa(v)^+$ is any element with $\bar v(x)=\gamma$,
it follows easily that the valuation topology on $\kappa(v)^+$
agrees with the $x$-adic topology, so $v$ is a Tate valuation.
Next, after replacing $A_0$ by the smallest $\Z$-subalgebra
of $A_0$ containing $I$, we may assume that
$\pi_v(A_0)\subset\kappa(v)^+$, and notice that
$\pi_v(I^n)$ generates $x^n\cdot\kappa(v)^+$ for every
$n\in\N$, so that $\pi_v$ is continuous and f-adic
(details left to the reader).

(iii): Say that $w=v^\Delta$ for some convex subgroup
$\Delta$ strictly contained in $\Gamma_{\!v}$. Then
the value $\gamma_v$ cannot lie in $\Delta$, so $w(I)=0$,
whence the assertion.

(iv): Proposition \ref{prop_reverse-order} says that $w$
is a primary specialization of a secondary specialization
$u$ of $v$, and $u$ is a continuous valuation of $A$, by
remark \ref{rem_Cont-A}(iv). In light of (iii), it follows
that $u=w$, whence the contention.

(v): First, from (iv) we see that every continuous generization
of $v$ must be a secondary generization and must be analytic;
hence, let $\Delta\subset\Gamma_{\!v}$ be the largest convex
subgroup that does not contain $\gamma_v$, and set $w:=v_\Delta$.
It is easily seen that the image $\bar\gamma_v$ of $\gamma_v$
in $\Gamma_{\!w}=\Gamma_{\!v}/\Delta$ is still final, and by
definition $\gamma_w=\bar\gamma_v$, so $w$ is continuous, by (i).
It remains only to check that $w$ does not admit any proper
(secondary) generization; however, by construction, every
non-trivial convex subgroup $\Delta'\subset\Gamma_{\!w}$ contains
$\gamma_w$, and on the other hand, $\gamma_w$ cannot lie in any
proper convex subgroup of $\Gamma_{\!w}$. Thus, the only convex
subgroups of $\Gamma_{\!w}$ are $\{1\}$ and $\Gamma_{\!w}$, whence
the assertion.
\end{proof}

\begin{theorem}\label{th_Cont-spectral}
Let $A$ be any f-adic ring, and $J\subset A$ any finitely
generated ideal such that $\Spec\,A/J$ is the non-analytic
locus of\/ $\Spec\,A$. Then we have :
\begin{enumerate}
\item
$\Cont(A)=\{v\in\Spv(A,J)~|~
         \text{$v(a)<1$ for every $a\in A^{\circ\circ}$}\}$.
\item
$\Cont(A)$ is a closed (in particular, pro-constructible)
subset of\/ $\Spv(A,J)$. Especially, $\Cont(A)$ is a spectral
topological space. 
\end{enumerate}
\end{theorem}
\begin{proof}(i): Let $A_0\subset A$ be a subring of
definition, and $I\subset A_0$ an ideal of adic definition;
by lemma \ref{lem_deja-vu}(i) we have $\Spec\,A/IA=\Spec\,A/J$,
so that the ideals $J$ and $IA$ have the same radical;
from lemma \ref{lem_was-rem-Spv-A-I}(i) we deduce that
$\Spv(A,J)=\Spv(A,IA)$, so we may assume from start
that $J=IA$. Since $I$ is an open subset of $A$, for
every $a\in A^{\circ\circ}$ we may find $n\in\N$ such
that $a^n\in I$; combining with lemma \ref{lem_Cont-A}(i)
we easily deduce that
$$
\Cont(A)=\{v\in\Spv\,A~|~
\text{$v(a)$ is final in $\Gamma_{\!v}$, for every
$a\in I$}\}.
$$
In light of lemma \ref{lem_cGamma-I}(i) we then see that
$$
\Cont(A)\subset\{v\in\Spv(A,J)~|~
\text{$v(a)<1$ for every $a\in I$}\}
$$
and to conclude, it suffices to prove the converse inclusion.
Hence, let $v\in\Spv(A,J)$ be any valuation such that $v(a)<1$
for every $a\in I$. If $\Gamma_{\!v}\neq c\Gamma_{\!v}$, lemma
\ref{lem_cGamma-I}(i) implies that either $v(a)=0$ or $v(a)$
is final in $\Gamma_{\!v}$ for every $a\in J$, and then
$v\in\Cont(A)$, by the foregoing.

Lastly, suppose that $\Gamma_{\!v}=c\Gamma_{\!v}$, and let $a\in I$
be any element; in this case, it suffices to show that, for
every $b\in A$ with $v(b)\neq 0$, there exists $n\in\N$ such
that $v(a)^n<v(b)^{-1}$. However, we may find $n\in\N$ such
that $a^nb\in I$, so $v(a^nb)<1$ by assumption, whence the
contention.

(ii): From (i) we see that $\Cont(A)=\bigcap_{a\in A^{\circ\circ}}
\bigl(\Spv(A,J)\setminus R_A\bigl(\frac{1}{a}\bigr)\bigr)$.
Combining with corollary \ref{cor_procon-is-spec}, we obtain
the assertion.
\end{proof}

\sset\subsubsection{}\label{subsec_rational-Cont}
Let $A$ be and $J$ be as in theorem \ref{th_Cont-spectral}.
We say that a subset of $\Cont(A)$ is {\em rational}, if it
can be written as an intersection $\Cont(A)\cap R$, where
$R$ is a rational subset of $\Spv(A,J)$ (see definition
\ref{def_Spv-I}(ii)). Hence, a rational subset is of the form
$\Cont(A)\cap R_A\bigl(\frac{f_1}{f_0},\dots,\frac{f_n}{f_0}\bigr)$,
where $(f_0,f_1,\dots,f_n)$ is any sequence of elements of
$A$ that generates an open ideal. Especially, the class of
rational subsets of $\Cont(A)$ does not depend on the choice
of the ideal $J$, and provides a basis of quasi-compact open
subsets of $\Cont(A)$ that is closed under finite intersections
(theorem \ref{th_Spv-I}(ii)).

\begin{corollary}\label{cor_Cont-spectral}
Let $f:A\to B$ be an f-adic morphism of f-adic topological
rings. We have:
\begin{enumerate}
\item
$\Cont(f)$ is a spectral map. More precisely, if
$R$ is a rational subset of $\Cont(A)$, then
$\Cont(f)^{-1}(R)$ is a rational subset of $\Cont(B)$.
\item
$\Cont(f)$ restricts to maps
$\Cont(B)_\mathrm{a}\to\Cont(A)_\mathrm{a}$ and
$\Cont(B)_\mathrm{na}\to\Cont(A)_\mathrm{na}$.
\end{enumerate}
\end{corollary}
\begin{proof}(i): Indeed, if
$R=\Cont(A)\cap R_A\bigl(\frac{a_1}{a_0},\dots,\frac{a_n}{a_0}\bigr)$
for a sequence $a_0,a_1,\dots,a_n$ that generates an open
ideal of $A$, then $\Cont(f)^{-1}(R)=\Cont(B)\cap
R_A\bigl(\frac{f(a_1)}{f(a_0)},\dots,\frac{f(a_n)}{f(a_0)}\bigr)$,
and since $f$ is adic, it is easily seen that the sequence
$f(a_0),f(a_1),\dots,f(a_n)$ generates an open ideal of $B$.

(ii) is clear.
\end{proof}

\begin{proposition}\label{prop_Cont-open-subring}
Let $A$ be an f-adic ring, $B\subset A$ an open subring,
and denote by $i:B\to A$ the inclusion map. Then we have :
\begin{enumerate}
\item
$\Cont(A)=\Spv(i)^{-1}\Cont(B)$.
\item
$\Cont(i)$ restricts to a homeomorphism
$\Cont(A)_\mathrm{a}\isom\Cont(B)_\mathrm{a}$.
\end{enumerate}
\end{proposition}
\begin{proof}(i): Let $v$ be a valuation of $A$ such
that $w:=\Spv(i)(v)=v_{|B}$ is a continuous valuation
of $B$; we have to show that $v$ is a continuous
valuation of $A$, and by remark \ref{rem_Cont-A}(vi),
we may assume that $\Ker\,v$ is not open. In this
case, lemma \ref{lem_deja-vu}(iii) implies that
$\kappa(v)=\kappa(w)$, and therefore $v$ and $w$
have the same value group; the assertion follows
immediately.

(ii): Define $X_A$ and $X_A^{\circ\circ}$ as in definition
\ref{def_deja-vu}; we have
$$
\Cont(A)_\mathrm{a}=\Spv(X_A\setminus X_A^{\circ\circ})\cap\Cont(A)
$$ and
likewise for $\Cont(B)_\mathrm{a}$, so the assertion follows
from (i) and lemma \ref{lem_deja-vu}(iii).
\end{proof}

\begin{proposition}\label{prop_Cont-of-complete}
Let $A$ be a topological ring, $A^\wedge$ the separated
completion of $A$. We have:
\begin{enumerate}
\item
The completion map $i:A\to A^\wedge$ induces a bijective
and continuous map
$$
\Cont(i):\Cont(A^\wedge)\to\Cont(A).
$$
\item
Moreover, with the ordering given by specializations
on $\Cont(A)$ and $\Cont(A^\wedge)$, the map $\Cont(i)$
is also an isomorphism of partially ordered sets
(see remark {\em\ref{rem_Cont-A}(v)})
$$
(\Cont(A^\wedge),\leq)\isom(\Cont(A),\leq).
$$
\item
Furthermore, $\Cont(i)$ restricts to bijections
$$
\Cont(A^\wedge)_\mathrm{a}\isom\Cont(A)_\mathrm{a}
\qquad\text{and}\qquad
\Cont(A^\wedge)_\mathrm{na}\isom\Cont(A)_\mathrm{na}.
$$
\item
If $A$ is f-adic, $\Cont(i)$ is a homeomorphism.
\end{enumerate}
\end{proposition}
\begin{proof}(i): For the injectivity of $\Cont(i)$
we show the following more precise :

\begin{claim}\label{cl_complete-cont-val}
Let $w\in\Cont(A^\wedge)$ be any element, and set
$v:=\Cont(i)(w)$. Let $\Gamma_{\!w}$ and $\Gamma_{\!v}$
be the value groups of $w$ and respectively $v$. Then we have:
\begin{enumerate}
\item
$\Gamma_{\!w}=\Gamma_{\!v}$ and $c\Gamma_{\!w}=c\Gamma_v$.
\item
The induced map of residue fields $\kappa(v)\to\kappa(w)$
has dense image, for the valuation topology of $\kappa(w)$
(notation of remark \ref{rem_semi-norm}(v)).
\end{enumerate}
\end{claim}
\begin{pfclaim}(i): Let $a\in A^\wedge$ be any element such
that $\gamma:=w(a)\in\Gamma_{\!w}$; by assumption, there
exists an open neighborhood $U_\gamma$ of $0$ in
$A^\wedge$ such that $w(x)<\gamma$ for every
$x\in U_\gamma$, hence $w(a+U_\gamma)=\{\gamma\}$,
and since $a+U_\gamma$ intersects the image of $A$ in
$A^\wedge$, we deduce that $\gamma\in\Gamma_{\!v}$, whence
the claim.

(ii): According to remark \ref{rem_Cont-A}(i), the natural
map $\pi_v:A\to\kappa(v)$ is continuous for the valuation
topology $\cT_v$ of $\kappa(v)$, hence it extends uniquely
to a continuous ring homomorphism
$\pi_v^\wedge:A^\wedge\to\kappa(v)^\wedge$, where
$(\kappa(v)^\wedge,\cT^\wedge_v)$ denotes the completion of
$(\kappa(v),\cT_v)$. However, $v$ extends uniquely to a
valuation $v^\wedge:\kappa(v)^\wedge\to\Gamma_{v\circ}$ and
$\cT^\wedge_v$ agrees with the corresponding valuation
topology of $\kappa(v)^\wedge$ (proposition
\ref{prop_stays-valuation}(iv,v)). It follows that
$v^\wedge\circ\pi_v^\wedge:A^\wedge\to\Gamma_{v\circ}$ is
the unique continuous valuation whose composition with
$i$ agrees with $v$. Consequently, there exists as well
a unique inclusion of fields $j:\kappa(w)\to\kappa(v)^\wedge$
whose composition with the projection
$\pi_w:A^\wedge\to\kappa(w)$ agrees with $\pi_v^\wedge$,
and $v^\wedge\circ j$ is the residual valuation $\cT_w$
of $\kappa(w)$. Especially, $\cT_w$ agrees with the topology
induced by $\cT^\wedge_v$ via $j$; since the image of
$\kappa(v)$ is dense in $\kappa(v)^\wedge$ relative to the
topology $\cT^\wedge_v$ (theorem \ref{th_complete-top-grps}(ii)),
the assertion follows.
\end{pfclaim}

Now, let $w,w'\in\Cont(A^\wedge)$ be two valuations such
that $w\circ i$ is equivalent to $w'\circ i$, and denote
by $\Gamma_{\!w}$ (resp. $\Gamma_{\!w'}$) the value group
of $w$ (resp. of $w'$). From claim \ref{cl_complete-cont-val}
we may assume that $\Gamma_{\!w}=\Gamma_{\!w'}$ and
$w\circ i=w'\circ i$. Since the image of $A$ is dense
in $A^\wedge$ and the topology of $\Gamma_{\!w}$ is separated,
we conclude that $w=w'$, as required.

Next, we show that $\Cont(i)$ is surjective. Indeed,
let $v$ be any continuous valuation of $A$ with value
group $\Gamma$, so that, for every $\gamma\in\Gamma$
there exists an open neighborhood $U_\gamma$ of $0$ in
$A$ such that $v(x)<\gamma$ for every $x\in U_\gamma$;
we need to show that there exists a valuation
$v^\wedge:A^\wedge\to\Gamma_{\!\circ}$ such that
$v=v^\wedge\circ i$. We define $v^\wedge$ as follows. Let
$\cC:=(a_i~|~i\in I)$ be any Cauchy net in $A$ (indexed
by some filtered ordered set $I$); for every
$\gamma\in\Gamma$, we may then find $i(\gamma)\in I$
such that $a_j-a_k\in U_\gamma$ for every $j,k\in I$ with
$j,k\geq i(\gamma)$. We let $\cC_\gamma$ be the Cauchy
net $(a_j~|~j\in I,\ j\geq i(\gamma))$ for every such
$\gamma$. Now, suppose that there exists $\gamma\in\Gamma$
and $a\in\cC_\gamma$ with $\delta:=v(a)\geq\gamma$; then
$v(b)=\delta$ for every $b\in\cC_\gamma$, and in this case
we set $v^\wedge(\cC):=\delta$. Otherwise, we set
$v^\wedge(\cC):=0$. It is easily seen that this definition
depends only on the equivalence class of $\cC$ in $A^\wedge$,
so it yields a well defined map
$v^\wedge:A^\wedge\to\Gamma_{\!\circ}$, and it is clear that
$v^\wedge\circ i=v$. We leave to the reader the verification
that $v^\wedge$ is indeed a continuous valuation on $A^\wedge$.

(ii): Let $u\in\Cont(A^\wedge)$ be any element, and $w'$
a specialization of $u':=\Cont(i)(u)$ in $\Cont(A)$. By
proposition \ref{prop_reverse-order} there exists
$v'\in\Spv\,A$ that is both a secondary specialization
of $u'$ and a primary generization of $w'$ in $\Spv\,A$.
By lemma \ref{lem_image-of-special}(iii) we may find
a secondary specialization $v$ of $u$ in $\Spv\,A^\wedge$
with $\Spv(i)(v)=v'$; then $v\in\Cont(A^\wedge)$ and
$v'\in\Cont(A)$ (remark \ref{rem_Cont-A}(iv)). By claim
\ref{cl_complete-cont-val}, we know that $v$ and $v'$
have the same value groups and the same characteristic
subgroups; hence, we may find a primary specialization
$w$ of $v$ in $\Spv\,A$ with $w'=\Spv(i)(w)$, and
$w\in\Cont(A)$ (remark \ref{rem_Cont-A}(iv)), as required.

(iii): We have already observed that $\Cont(i)$ restricts
to a map $\Cont(A^\wedge)_\mathrm{na}\to\Cont(A)_\mathrm{na}$
(remark \ref{rem_Cont-A}(vii)); in light of (i), it then
suffices to show that $\Cont(i)$ restricts as well to a map
$\Cont(A^\wedge)_\mathrm{a}\to\Cont(A)_\mathrm{a}$. Thus, let
$v:A^\wedge\to\Gamma_\circ$ be any continuous analytic valuation,
and suppose that $\fp:=\Ker\,v\circ i$ is open in $A$; by
continuity, the topological closure $i(\fp)^c$ of $i(\fp)$
in $A^\wedge$ lies in $\Ker\,v$. However, $i(\fp)^c$ contains an
open subset of $A^\wedge$ (claim \ref{cl_maximality-of-U-wedge}(ii)),
and therefore it is open, a contradiction.

(iv): From (i), we know already that $\Cont(i)$ is continuous
and bijective. It then suffices to show that $\Cont(i)$ is a
closed map, if $A$ is f-adic. However, in this case $\Cont(A)$
and $\Cont(A^\wedge)$ are spectral topological spaces (theorem
\ref{th_Cont-spectral}(ii)) and $\Cont(i)$ is a spectral map
(corollary \ref{cor_Cont-spectral}(i) and proposition
\ref{prop_complete-f-adic}(ii)), so it suffices to check
that $\Cont(i)$ is specializing (proposition
\ref{prop_closed-under-spec}(i) and corollary
\ref{cor-pro-constr}(i)). The latter follows immediately
from (ii).
\end{proof}

\begin{corollary}\label{cor_characterize-f-adic}
Let $f:A\to B$ be a continuous ring homomorphism of
f-adic rings. Suppose that $B$ is topologically local,
and denote by $B^\wedge$ the separated completion of $B$.
Then the following conditions are equivalent :
\begin{enumerate}
\alphaenu
\item
$f$ is an f-adic map.
\item
The composition $A\to B^\wedge$ of $f$ with the completion
map $B\to B^\wedge$ is f-adic.
\item
$\Cont(f)$ restricts to a map
$\Cont(B)_\mathrm{a}\to\Cont(A)_\mathrm{a}$.
\end{enumerate}
\end{corollary}
\begin{proof}(a)$\Rightarrow$(b) since the completion
map $B\to B^\wedge$ is f-adic (proposition
\ref{prop_complete-f-adic}(i,ii)).

(b)$\Rightarrow$(c) by corollary \ref{cor_Cont-spectral}(ii)
and proposition \ref{prop_Cont-of-complete}(i).

(c)$\Rightarrow$(a): Let $B_0$ be a subring of definition
of $B$, and $A_0$ a subring of definition of $A$ such that
$A_0\subset f^{-1}B_0$ (corollary \ref{cor_f-adics}(ii)),
and let $f_0:A_0\to B_0$ be the restriction of $f$.
From (c) and proposition \ref{prop_Cont-open-subring}(ii)
we deduce that $\Cont(f_0)$ restricts to a map
$\Cont(B_0)_\mathrm{a}\to\Cont(A_0)_\mathrm{a}$, and it
suffices to check that $f_0$ is f-adic. Moreover,
$B_0$ is still topologically local (proposition
\ref{prop_quasi-affinoid}(i)). Thus, we may replace
$f$ by $f_0$, and assume from start that $B=B^\circ$,
so that $B^{\circ\circ}$ lies in the Jacobson radical
of $B$. Now, define $X_A$ and $X^{\circ\circ}_A$ as in
definition \ref{def_deja-vu}, and likewise for the
ring $B$; also, set $\phi:=\Spec\,f$. By lemma
\ref{lem_deja-vu}(iv), it suffices to show that
$\phi^{-1}X_A^{\circ\circ}=X^{\circ\circ}_B$. However, suppose
that the latter fails; then there exists a prime ideal
$\fp\in X_B\setminus X^{\circ\circ}_B$ with
$\phi(\fp)\in X^{\circ\circ}_A$, and we pick any maximal
ideal $\fm$ of $B$ containing $\fp$. Set $C:=(B/\fp)_\fm$,
and choose a valuation ring $V$ of $\Frac\,C$ that
dominates $C$ (corollary \ref{cor_cornerstone}). The
valuation ring $V$ corresponds to a point $v\in\Spv\,B$;
next, pick any finitely generated ideal of definition
$I$ of $B$, and set $w:=v^{c\Gamma_{\!v}(I)}$. Since
$B^{\circ\circ}\subset\fm$, theorem \ref{th_Cont-spectral}(i)
implies that $w$ is a continuous valuation of $B$.
Since $\fp=\Ker\,v$ and $I$ is not contained in $\fp$,
it follows that $I$ is not contained in the support
of $w$, therefore $w\in\Cont(B)_\mathrm{a}$. However,
$w':=\Cont(f)(w)$ is a primary specialization of
$\Spv\,f(v)$, and the support of the latter is an open
ideal of $A$, by construction, so $w'\in\Cont(A)_\mathrm{na}$,
which contradicts our assumptions.
\end{proof}

\begin{proposition} Let $A$ be any f-adic topological ring,
and denote by $A^\sep$ the maximal separated quotient of $A$.
The following holds :
\begin{enumerate}
\item
$\Spv\,A=\emptyset$ if and only if $A=0$.
\item 
$\Cont(A)=\emptyset$ if and only if $A^\sep=0$.
\item
$\Cont(A)_\mathrm{a}=\emptyset$ if and only if the topology
of $A^\sep$ is discrete.
\end{enumerate}
\end{proposition}
\begin{proof}(i): Indeed, if $A\neq 0$, it has a
prime ideal $\fp$, and the trivial valuation on
$A/\fp$ yields an element of $\Spv\,A$.

(iii): Let $A^\wedge$ be the separated completion of $A$;
if $\Cont(A)_\mathrm{a}=\emptyset$, we have
$\Cont(A^\wedge)_\mathrm{a}=\emptyset$ as well, by
proposition \ref{prop_Cont-of-complete}(iii), and it
suffices to show that the topology of $A^\wedge$
is discrete. Let $B\subset A^\wedge$ be a subring of
definition, and $I\subset B$ a finitely generated
ideal of adic definition (proposition
\ref{prop_complete-f-adic}(i)); then
\set\begin{equation}\label{eq_Cont-is-empty}
\Cont(B)_\mathrm{a}=\emptyset
\end{equation}
as well (proposition \ref{prop_Cont-open-subring}(ii)),
and we remark :

\begin{claim}\label{cl_I-adic-discrete}
$I$ is a nilpotent ideal. 
\end{claim}
\begin{pfclaim} Let $\fp\subset B$ be any prime ideal,
and suppose by way of contradiction, that $I$ is not
contained in $\fp$. Choose a maximal ideal $\fm$ of
$C:=B/\fp$, and pick a valuation ring $V$ of
$\Frac\,C$ that dominates $C_\fm$ (corollary
\ref{cor_cornerstone}).
Let $\fn\subset B$ be the maximal ideal such that
$\fn/\fp=\fm$; then $V$ corresponds to a valuation
$v$ of $B$ with support equal to $\fp$, and such that
$v(a)<1$ for every $a\in\fn$. Since $B$ is $I$-adically
complete, we have $I\subset\fn$ (remark
\ref{rem_someth-on-bdd-in-Z-lin}(v)),
so it follows easily that $v(a)<1$ for every
$a\in B^{\circ\circ}$. Let $\Gamma_{\!v}$ be the value group
of $v$, and set $w:=v^{c\Gamma_{\!v}(I)}$ (notation of
definition \ref{def_cGamma-I}). Then $w\in\Spv(B,I)$
and $w(a)<1$ for every $a\in B^{\circ\circ}$, so $w$ is
continuous (theorem \ref{th_Cont-spectral}(i)). Lastly,
since $v(I)\neq\{0\}$ by construction, we have $w(I)\neq\{0\}$
as well, since $w$ is an $I$-admissible specialization
of $v$. Therefore $w\in\Cont(B)_\mathrm{a}$, which
contradicts \eqref{eq_Cont-is-empty}.
\end{pfclaim}

Claim \ref{cl_I-adic-discrete} means that the topology
of $B$ is discrete, hence the same holds for the topology
of $A^\wedge$, as required.

(ii): From (iii) we know that the topological closure
$\{0\}^c$ of the ideal $\{0\}$ of $A$ is open in $A$;
hence, every prime ideal containing $\{0\}^c$ is open
as well, and therefore every valuation with support
equal to such a prime ideal would be continuous
(remark \ref{rem_Cont-A}(vi)); since $\Cont(A)=\emptyset$,
we conclude that no prime ideal of $A$ contains
$\{0\}^c$, {\em i.e.} $\{0\}^c=A$, whence the assertion.
\end{proof}

\begin{lemma}\label{lem_criterion-opennes}
Let $A$ be a topologically local f-adic ring, $J\subset A$
an ideal. The following conditions are equivalent :
\begin{enumerate}
\alphaenu
\item
$J$ is an open ideal.
\item
For every $v\in\Cont(A)_\mathrm{a}$ we have $J\not\subset\Ker\,v$.
\item
For every $v\in\Cont(A)_\mathrm{a}$ of rank one, we have
$J\not\subset\Ker\,v$.
\end{enumerate}
\end{lemma}
\begin{proof} Clearly (a)$\Rightarrow$(b)$\Rightarrow$(c).

Next, suppose that $J$ is not open; to conclude, it suffices
to exhibit a rank one continuous analytic valuation on $A$ whose
support contains $J$. Now, let $B\subset A$ be a subring of
definition, and set $J_B:=B\cap J$; since $J$ is not open
in $A$, the ideal $J_B$ cannot be open in $B$. Moreover,
$B$ is also topologically local, by proposition
\ref{prop_quasi-affinoid}(i).

\begin{claim}\label{cl_reduce-to-B}
It suffices to exhibit a rank one continuous analytic
valuation on $B$ whose support contains $J_B$. 
\end{claim}
\begin{pfclaim} Indeed, let $v:B\to\Gamma_\circ$ be such a
valuation, and denote by $\fp\subset B$ the support of $v$;
with the notation of definition \ref{def_deja-vu}, we have
$\fp\in X_B\setminus X^{\circ\circ}_B$, and by lemma
\ref{lem_deja-vu}(iii) the inclusion map $i:B\to A$ induces
an isomorphism of schemes
$\phi:X_A\setminus X^{\circ\circ}_A\isom X_B\setminus X^{\circ\circ}_B$.
Set $\fq:=\phi^{-1}(\fp)$; it follows that $i$ extends to a
ring isomorphism $i_\fp:B_\fp\isom A_\fq$. Let also $j:A\to A_\fq$
be the localization map; $v$ extends uniquely to a valuation
$\bar v:B_\fp\to\Gamma_\circ$, and composing with $i_\fp^{-1}\circ j$
we deduce a valuation $w:A\to\Gamma_\circ$. Moreover,
$\Cont(i):\Cont(A)\to\Cont(B)$ is a homeomorphism
(proposition \ref{prop_Cont-open-subring}(ii)), and a direct
inspection of the definitions yields
$$
\Cont(i)^{-1}(v)=w.
$$
Especially, $w$ is continuous, analytic and of rank one.
Lastly, the support of $\bar v$ contains
$J_{B,\fp}=i_\fp^{-1}J_\fq$, so the support of $w$ contains $J$.
\end{pfclaim}

Due to claim \ref{cl_reduce-to-B}, we may replace $A$ by
$B$ and $J$ by $J_B$, and assume from start that $A=A^\circ$;
especially, $A^{\circ\circ}$ lies in the Jacobson radical of
$A$. By lemma \ref{lem_deja-vu}(v), there exists
$x\in\Spec\,A/J\setminus X^{\circ\circ}_A$; according to
remark \ref{rem_specialize}(ii) we may find a minimal
specialization $y$ of $x$ in
$\Spec\,A/J\setminus X^{\circ\circ}_A$, and if $\{y\}^c$
denotes the topological closure of $\{y\}$ in
$\Spec\,A/J$, the closed subset
$Z:=\{y\}^c\cap X^{\circ\circ}_A$ of $\Spec\,A/J$ is
non-empty, since $X_A^{\circ\circ}$ contains all the
maximal ideals of $A$. We let $\fp\subset A$ be any
maximal point of $Z$, and we notice that the image
$C$ of $A_\fp$ in $\kappa(y)$ is a one-dimensional
local domain.

\begin{claim}\label{cl_one-dim-dominates}
For any one-dimensional local domain $(R,\fm_R)$
there exists a one-dimensional valuation ring of
$\Frac\,R$ that dominates $R$.
\end{claim}
\begin{pfclaim} By corollary \ref{cor_cornerstone} we
may find a valuation ring $V$ of $\Frac\,R$ that
dominates $R$. Let $\cF$ be the set of all prime
ideals $\fq$ of $V$ such that $\fq\cap R=\fm_R$;
then $\cF$ is non-empty, and it admits a minimal
element $\fr:=\bigcap_{\fq\in\cF}\fq$. It then suffices to
check that the localization $V_\fr$ is one-dimensional.
However, let $\fp\subset V_\fr$ be any non-zero
prime ideal, and pick any $a\in\fp\setminus\{0\}$;
we may write $a:=b/b'$ for some $b,b'\in R$, and
then $b\in\fp$ as well, so $\fp\cap R=\fm_R$, in
which case $\fp$ must be the maximal ideal of $V_\fr$,
and the assertion follows.
\end{pfclaim}

By claim \ref{cl_one-dim-dominates}, we may find a
one-dimensional valuation ring $V$ of $\kappa(y)$
that dominates $C$, and denote by $v\in\Spv\,A$ the
corresponding valuation. Let also $I\subset A$ be
any finitely generated ideal of adic definition; by
construction, the support of $v$ contains $J$ and is
not an open ideal, and $v$ does not admit any proper
$I$-admissible specializations (see
\eqref{subsec_I-admissible}), so $v\in\Spv(A,I)$.
Moreover $A^{\circ\circ}\subset\fp$, so that $v(b)<1$
for every $b\in A^{\circ\circ}$, and consequently
$v\in\Cont(A)_\mathrm{a}$. 
\end{proof}

\subsection{Affinoid rings and affinoid schemes}
In this section we present an extension and refinement of
Huber's theory of affinoid rings.

\begin{definition}\label{def_affinoid}
(i)\ \
Let $A$ be any f-adic ring. An open subring $B\subset A$ is
called a {\em ring of integral elements of $A$}, if
$B\subset A^\circ$ and $B$ is integrally closed in $A$.

(ii)\ \
An {\em affinoid ring} is a datum $(A,A^+)$ consisting of
an f-adic ring $A$ and a ring $A^+$ of integral elements
of $A$. A {\em morphism of affinoid rings} $f:(A,A^+)\to(B,B^+)$
is a continuous ring homomorphism $f:A\to B$ such that
$f(A^+)\subset B^+$. Then we say that $f$ is {\em f-adic},
if the same holds for the underlying map $A\to B$.

(iii)\ \
A {\em quasi-affinoid ring} is a datum $(A,A^+,U)$
consisting of an affinoid ring $(A,A^+)$ and a constructible
open subset $U\subset X_A$ that contains the analytic locus
(see definition \ref{def_deja-vu}).
A {\em morphism of quasi-affinoid rings} $f:(A,A^+,U)\to(B,B^+,V)$
is a morphism $f:(A,A^+)\to(B,B^+)$ of affinoid rings
such that $\Spec\,f$ restricts to a morphism of
schemes $V\to U$. We say that such $f$ is {\em f-adic},
if the same holds for the underlying map $A\to B$.

(iv)\ \
A {\em quasi-affinoid scheme} is a datum
$\underline X:=(X,\cT_X,A_X^+)$ consisting of a (quasi-compact)
quasi-affine scheme $X$ (see \cite[Ch.II, D\'ef.5.1.1]{EGAII}),
an f-adic topology $\cT_X$ on $A_X:=\cO_X(X)$, and a subring
$A_X^+\subset A_X$ such that the following holds :
\begin{itemize}
\item
$(A_X,A_X^+)$ is an affinoid ring.
\item
For every $s\in A_X^{\circ\circ}$, the open subscheme
$X_s:=X\times_{\Spec\,A_X}\Spec\,A_X[s^{-1}]$ is affine.
\end{itemize}
A {\em morphism of quasi-affinoid schemes}
$\phi:(X,\cT_X,A_X^+)\to(Y,\cT_Y,A_Y^+)$ is a morphism
of schemes $\phi:X\to Y$ whose associated morphism of
sheaves of rings $\phi^\flat:\cO_Y\to\phi_*\cO_X$ induces
a morphism of affinoid rings
$\phi^\flat_Y:(\cO_Y(Y),A_Y^+)\to(\cO_X(X),A_X^+)$.
We say that $\phi$ is {\em f-adic} if $\phi^\flat_Y$
is f-adic. We say that $\underline X$ is an
{\em affinoid scheme} if $X$ is an affine scheme.

(v)\ \
Let $\underline A:=(A,A^+,U)$ (resp.
$\underline X:=(X,\cT_X,A^+_X)$) be any quasi-affinoid ring
(resp. quasi-affinoid scheme). We say that $\underline A$
(resp. $\underline X$) is {\em topologically local} (resp.
{\em topologically henselian}, resp. {\em complete},
resp. {\em separated}) if the underlying f-adic ring $A$
(resp. $\cO_{\!X}(X)$) enjoys the corresponding property.
\end{definition}

\begin{remark}\label{rem_shouldbe-quasi-affinoid}
(i)\ \
Clearly the affinoid rings, the quasi-affinoid rings,
the affinoid schemes and the quasi-affinoid schemes with
their morphisms as in definition \ref{def_affinoid} form
categories denoted respectively : 
$$
\mathsf{Afd.Ring}
\qquad
\mathsf{q.Afd.Ring}
\qquad
\mathsf{Afd.Sch}
\qquad
\mathsf{q.Afd.Sch}
$$
($\mathsf{Afd.Sch}$ is a full subcategory of $\mathsf{q.Afd.Sch}$).
Also, we have a natural fully faithful functor
$$
\mathsf{Afd.Ring}\to\mathsf{q.Afd.Ring}
\qquad
(A,A^+)\mapsto(A,A^+,\Spec\,A)
$$
so we may regard $\mathsf{Afd.Ring}$ as a full subcategory
of $\mathsf{q.Afd.Ring}$.

(ii)\ \
By virtue of \cite[Ch.II, Prop.5.1.2]{EGAII}, we have a
natural fully faithful functor
$$
\sGamma:\mathsf{q.Afd.Sch}^o\to\mathsf{q.Afd.Ring}
\qquad
(X,\cT_X,A_X^+)\mapsto(A_X,A_X^+,X)
$$
(where $A_X:=\cO_X(X)$ and $X$ is naturally identified
with the image of the open immersion $X\to\Spec\,A_X$).
Indeed, let $s\in A_X^{\circ\circ}$ be any element, and
define the open subscheme $X_s\subset X$ as in definition
\ref{def_affinoid}(iv); by assumption, $X_s$ is affine,
and on the other hand we have a natural isomorphism
of $A_X$-algebras
$$
\Gamma(X_s,\cO_X)\isom A_X[s^{-1}]
$$
by \cite[Ch.I, Cor.9.2.2]{EGAI}. Therefore the image
of the open immersion $\Spec\,A_X[s^{-1}]\to\Spec\,A_X$
equals $X_s$; especially, $\Spec\,A_X[s^{-1}]\subset X$
for every $s\in A_X^{\circ\circ}$, so $(A_X,A_X^+,X)$ is
indeed a quasi-affinoid ring. Especially, the analytic
locus of $\Spec\,A_X$ lies in $X$, and shall also be called
the {\em analytic locus of $\underline X$}.

(iii)\ \
Let $(A,A^+)$ be any affinoid ring, $B\subset A$ an
open subring, and endow $B$ with the topology induced
by the inclusion map $B\to A$. Then corollary
\ref{cor_Tate}(i) implies that the pair $(B,A^+\cap B)$
is an affinoid ring.

(iv)\ \
Let $A$ be any f-adic ring, and endow
$R:=\Z\oplus A^{\circ\circ}$ with the multiplication
map given by the rule : $(n,a)\cdot(m,b):=(nm,ma+nb+ab)$
for every $n,m\in\Z$ and $a,b\in A^{\circ\circ}$. Then $R$
is a ring, and we have a unique ring homomorphism
$f:R\to A$ that restricts to the inclusion map on
$A^{\circ\circ}$. The integral closure $B$ of the image of
$f$ is the smallest ring of integral elements of $A$.
On the other hand, $A^\circ$ is the largest subring of
integral elements of $A$.
\end{remark}

\begin{proposition}\label{prop_adjoint-to-sGamma}
{\em(i)}\ \
The functor $\sGamma$ of remark
{\em\ref{rem_shouldbe-quasi-affinoid}(ii)} admits a left adjoint :
$$
\sSpec:\mathsf{q.Afd.Ring}\to\mathsf{q.Afd.Sch}^o
\qquad
(A,A^+,U)\mapsto(U,\cT_U,A_U^+)
$$

{\em(ii)}\ \
If $f:\underline A\to\underline B$ is an f-adic morphism
of quasi-affinoid rings, then
$\sSpec(f):\sSpec(\underline B)\to\sSpec(\underline A)$
is an f-adic morphism of quasi-affinoid rings.

{\em(iii)}\ \
The functor $\sSpec$ is not fully faithful, but it
restricts to an equivalence of categories
$$
\sSpec:\mathsf{Afd.Ring}\to\mathsf{Afd.Sch}^o.
$$
\end{proposition}
\begin{proof}(i): For a quasi-affinoid ring
$\underline A:=(A,A^+,U)$, we let $\cT_U$ be the f-adic
topology on $A_U:=\cO_U(U)$ provided by proposition
\ref{prop_top-on-opens-fadic-case}(i), and let $A_U^+$ be
the integral closure of $A^+$ in $A_U$. Indeed, notice that
the restriction map $A\to A_U$ sends $A^\circ$ into $A^\circ_U$
(lemma \ref{lem_f-adics}(iii.a)); since $A^+\subset A^\circ$,
remark \ref{rem_someth-on-bdd-in-Z-lin}(iii) easily implies
that $A^+_U\subset A_U^\circ$. Next, let $\bar A:=\Img(A\to A_U)$
and endow $\bar A$ with the topology induced by $A_U$ via
the open inclusion map $\bar A\to A_U$; since $A^{\circ\circ}$
is an open subgroup of $A$, for every $s\in A_U^{\circ\circ}$
we may find $n\in\N$ and $t\in A^{\circ\circ}$ such that $s^n$
equals the image of $t$ in $\bar A$. Notice that $t$
annihilates the kernel of the projection $\pi:A\to\bar A$,
so $\Spec\,\pi$ restricts to an isomorphism of schemes
$\Spec\,\bar A[s^{-n}]\isom\Spec\,A[t^{-1}]\subset U$.
It follows that
$$
U\times_{\Spec\,A_U}\Spec\,A_U[s^{-1}]=U\times_{\Spec\,A}\Spec\,A[t^{-1}]
=\Spec\,A[t^{-1}]
$$
so the construction of $\sSpec(\underline A)$ is achieved.
Next, let $f:\underline A\to\underline B:=(B,B^+,V)$ be any
morphism of quasi-affinoid rings; set
$(V,\cT_V,B^+_V):=\sSpec(\underline B)$ and denote by
$\overline B$ the image of $B$ in $B_V:=\cO_V(V)$.
We deduce easily from proposition
\ref{prop_top-on-opens-fadic-case}(i) that (a) the subrings
$\overline A$ and $\overline B$ are open in $A_U$ and $B_V$,
so the induced map $\phi:A_U\to B_V$ is continuous if and only
the same holds for its restriction
$\overline f:\overline A\to\overline B$, and (b) the topologies
of $\overline A$ and $\overline B$ are induced from the surjections
$A\to\overline A$ and $B\to\overline B$, so the continuity of
$\overline f$ follows from that of $f$. Thus, $\phi$ is
continuous, and clearly $\phi(A^+_U)\subset B^+_V$, so
$\phi:\sSpec(\underline B)\to\sSpec(\underline A)$ is a well defined
morphism of quasi-affinoid schemes, and we let $\sSpec(f):=\phi$.

Now, a simple inspection shows that $\sSpec\circ\sGamma$
is the identity automorphism of the category
$\mathsf{q.Afd.Sch}^o$; on the other hand, for every
quasi-affinoid ring $\underline A:=(A,A^+,U)$, the
restriction map $A\to A_U$ determines a natural morphism
$\eta_{\underline A}:\underline A\to\sGamma\circ\sSpec(\underline A)$.
Lastly, it is easily seen that the pair of natural
transformations $(\eps:=\one_{\sSpec\circ\sGamma},\eta)$ fulfills the
triangular identities of \eqref{subsec_adj-pair}, so the assertion
follows from proposition \ref{prop_triangular-identities}(i).

(ii): Indeed, set $A_U:=\cO_U(U)$ and $B_V:=\cO_V(V)$; by assumption
there exist open and bounded subrings $A_0\subset A$
and $B_0\subset B$ such that $f(A_0)\subset B_0$, and
$f$ restricts to an adic ring homomorphism $A_0\to B_0$.
Then, let $\overline A_0$ be the image of $A_0$ in $A_U$,
and define likewise $\overline B_0\subset B_V$; we know
that $\overline A_0$ and $\overline B_0$ are open and
bounded subrings of $A_U$ and respectively $B_V$, for
the topologies induced by the latter topological rings,
and the image of every ideal of adic definition of $A_0$
(resp. of $B_0$) is an ideal of adic definition of $\overline A_0$
(resp. of $\overline B_0$) (proposition
\ref{prop_top-on-opens-fadic-case}(i)). It follows easily
that the continuous map $A_U\to B_V$ restricts to an adic
ring homomorphism $\overline A_0\to\overline B_0$, whence
the contention.

(iii) follows easily from the definitions.
\end{proof}

\begin{example}\label{ex_restriction-of-qaff}
Let $\underline X:=(X,\cT_X,A^+_X)$ be any quasi-affinoid
scheme, set $A_X:=\cO_{\!X}(X)$, and let $U\subset X$ be
any quasi-compact open subset containing the analytic
locus of $\Spec\,A_X$. Let us denote by $F_U$ the presheaf
on the category $\mathsf{q.Afd.Sch}$ such that
$$
F_U(\underline Z):=
\{\phi:\underline Z\to\underline X~|~\phi(Z)\subset U\}
\qquad
\text{for every quasi-affinoid scheme
$\underline Z:=(Z,\cT_Z,A^+_Z)$}
$$
so that $F_U$ is a sub-presheaf of the presheaf $h_{\underline X}$
represented by $\underline X$. Let $\underline A$ be the
quasi-affinoid ring $(A_X,A^+_X,U)$; we claim that $F_U$ is
representable by the quasi-affinoid scheme
$$
U\times_X\underline X:=\sSpec\,\underline A.
$$
Indeed, the identity map of $A_X$ yields a morphism of
quasi-affinoid rings $i:\sGamma(\underline X)\to\underline A$;
now, for every quasi-affinoid scheme
$\underline Z:=(Z,\cT_Z,A^+_Z)$, the datum of a morphism
$\psi:\underline Z\to U\times_X\underline X$ of quasi-affinoid
schemes is equivalent, by adjunction, to that of a morphism
$g:\underline A\to\sGamma(\underline Z)$ of quasi-affinoid
rings, and the composition
$f:=g\circ i:\sGamma(\underline X)\to\sGamma(\underline Z)$
is a morphism such that $\Spec(f)(Z)\subset U$, whence the
morphism $\phi:=\sSpec(f):\underline Z\to\underline X$ with
$\phi(Z)\subset U$. Conversely, any morphism
$\phi:\underline Z\to\underline X$ with $\phi(Z)\subset U$
yields a morphism
$\sGamma(\phi):\sGamma(\underline X)\to\sGamma(\underline Z)$
that factors uniquely through $\underline A$, whence the
contention.
\end{example}

\begin{example}\label{ex_f-adic-push-out}
(i)\ \
Consider two f-adic morphisms of affinoid rings
$$
\underline B':=(B',B'^+)\xleftarrow{\ g\ }
\underline B:=(B,B^+)\xrightarrow{\ f\ }
\underline A:=(A,A^+)
$$
and endow $A':=B'\otimes_BA$ with the f-adic topology $\cT_{A'}$
characterized by proposition \ref{prop_f-adic-push-out}(i),
so that the natural maps $f':B'\to A'$ and $g_A:A\to A'$
are both f-adic. Let also $R\subset A'$ be the image of
$B'^+\otimes_{B^+}A^+$, and denote by $A'^+$ the integral
closure of $R$ in $A'$. Then the datum
$$
\underline B'\otimes_{\underline B}\underline A:=(A',A'^+)
$$
is an affinoid ring. Indeed, a simple inspection of the
definition of $\cT_{A'}$ shows that $R$ is open in $A'$,
and we have $g_A(A^+),f(B^+)\subset A'^\circ$ by lemma
\ref{lem_f-adics}(iii), so also $R\subset A'^\circ$, and
then $A'^+\subset A'^\circ$ as well, by remark
\ref{rem_someth-on-bdd-in-Z-lin}(iv).

(ii)\ \
In the situation of (i), let also $U_B\subset X_B:=\Spec\,B$,
$U_A\subset X_A:=\Spec\,A$ and $U_{B'}\subset X_{B'}:=\Spec\,B'$
be three open subsets, such that $f$ and $g$ are morphisms
of quasi-affinoid rings
$$
(\underline B',U_{B'})\xleftarrow{\ g\ }
(\underline B,U_B)\xrightarrow{\ f\ }
(\underline A,U_A)
$$
and set $U_{A'}:=U_{B'}\times_{U_B}U_A$. Then the datum
$$
(\underline B',U_{B'})\otimes_{(\underline B,U_B)}(\underline A,U_A):=
(\underline B'\otimes_{\underline B}\underline A,U_{A'})
$$
is a quasi-affinoid ring as well.
Indeed, let $B_0\subset B$ be a ring of definition, and
$I_0\subset B_0$ an ideal of adic definition; we need to
show that $\Spec\,A'\setminus\Spec\,A'/I_0A'\subset U_{A'}$.
Set $Z_{B'}:=X_{B'}\setminus U_{B'}$ and $Z_A:=X_A\setminus U_A$;
the assertion is equivalent to
$(Z_{B'}\times_{X_B}X_A)\cup(X_{B'}\times_{X_B}Z_A)\subset
\Spec\,A'/I_0A'$. The latter is clear, since by assumption
$Z_{B'}\subset\Spec\,B'/I_0B'$ and $Z_A\subset\Spec\,A/I_0A$.
\end{example}

\begin{lemma}\label{lem_f-adic-pushout}
With the notation of example {\em\ref{ex_f-adic-push-out}},
the resulting commutative diagrams
$$
\xymatrix{ \underline B \ar[r]^-f \ar[d]_g &
\underline A \ar[d]^{g_A} &
(\underline B,U_B) \ar[r]^-f \ar[d]_g &
(\underline A,U_A) \ar[d]^{g_A} \\
\underline B' \ar[r]^-{f'} &
\underline  B'\otimes_{\underline B}\underline A &
(\underline B',U_{B'}) \ar[r]^-{f'} &
(\underline B',U_{B'})\otimes_{(\underline B,U_B)}(\underline A,U_A)
}$$
are cocartesian in the categories $\mathsf{Afd.Ring}$ and
respectively $\mathsf{q.Afd.Ring}$.
\end{lemma}
\begin{proof} Indeed, consider any two morphisms
$h:(A,A^+)\to(C,C^+)$, $k:(B,B^+)\to(C,C^+)$ of affinoid
rings such that $h\circ f=k\circ g$; by proposition
\ref{prop_f-adic-push-out}(i) there exists a unique
continuous ring homomorphism $l:A'\to C$ such that
$l\circ g_A=h$ and $l\circ f'=k$, and it is easily seen
that $l(A'^+)\subset C^+$, {\em i.e.} $l:(A',A'^+)\to(C,C^+)$
is a morphism of affinoid rings, whence the assertion
for the left diagram. The assertion for the right diagram
is an immediate consequence.
\end{proof}

\begin{example}\label{ex_shouldbe-quasi-affinoid}
Consider a quasi-affinoid ring $\underline A:=(A,A^+)$ and
any two f-adic morphisms $f:B\to A$, $g:B\to B'$ of f-adic
rings.

(i)\ \
As a special case of example \ref{ex_f-adic-push-out}(i),
take for $B^+$ and $B'^+$ the smallest subrings of integral
elements of $B$ and respectively $B'$ (see remark
\ref{rem_shouldbe-quasi-affinoid}(iv)) and set
$\underline B:=(B,B^+)$ and $\underline B':=(B',B'^+)$; we
may regard $f$ and $g$ as morphisms $\underline B\to\underline A$
and $\underline B\to\underline B'$ respectively. The resulting
tensor product $\underline B'\otimes_{\underline B}\underline A$
shall be simply denoted
$$
B'\otimes_B\underline A.
$$
By lemma \ref{lem_f-adic-pushout}, the resulting morphisms
$\underline A\xrightarrow{g_A}B'\otimes_B\underline A
\xleftarrow{f'}\underline B'$ enjoy the following universal
property. Let $(h,k)$ be any pair consisting of a morphism
$h:\underline A\to(C,C^+)$ of affinoid rings, and a continuous
ring homomorphism  $k:B'\to C$ such that $h\circ f=k\circ g$;
then there exists a unique morphism
$l:B'\otimes_B\underline A\to(C,C^+)$ of affinoid rings
such that $l\circ g_A=h$ and $l\circ f'=k$.

(ii)\ \
Likewise, let $(\underline A,U)$ be a quasi-affinoid ring.
As a special case of example \ref{ex_f-adic-push-out}(ii),
we may define $\underline B$ and $\underline B'$ as in (i),
and take $U_B:=\Spec\,B$, $U_{B'}:=\Spec\,B'$. There follow
morphisms of quasi-affinoid rings
$f:(\underline B,U_B)\to(\underline A,U)$ and
$f:(\underline B,U_B)\to(\underline B',U_{B'})$. The resulting
tensor product $(B'\otimes_B\underline A,U_{B'})$ will be
denoted simply
$$
B'\otimes_B(\underline A,U)
$$
and notice that $U_{B'}=\Spec\,A'\times_{\Spec\,A}U$.
It follows easily from lemma \ref{lem_deja-vu}(iv) that
this quasi-affinoid ring enjoys the corresponding universal
property as in (i).
Namely, let $(h,k)$ be any pair consisting of a morphism
$h:(\underline A,U)\to(C,C^+,V)$ of quasi-affinoid rings,
and $k:B'\to C$ a continuous ring homomorphism such that
$h\circ f=k\circ g$; then there exists a unique morphism
$l:B'\otimes_B(\underline A,U)\to(C,C^+,V)$ of quasi-affinoid
rings such that $l\circ g_A=h$ and $l\circ f'=k$.
\end{example}

\begin{example}\label{ex_fibre-prod-in-qaff-sch}
Let $\underline Y\to\underline X$ and
$\underline Y'\to\underline X$ be two f-adic morphisms
of quasi-affinoid schemes. Then the fibre product
$\underline Y\times_{\underline X}\underline Y'$ is
representable in the category $\mathsf{q.Afd.Sch}$.
Indeed, the induced morphisms of quasi-affinoid rings
$\sGamma(\underline X)\to\sGamma(\underline Y)$ and
$\sGamma(\underline X)\to\sGamma(\underline Y')$ are
f-adic, hence we may form the tensor product
$\underline A:=\sGamma(\underline Y)\otimes_{\sGamma(\underline X)}
\sGamma(\underline Y')$ as in example
\ref{ex_f-adic-push-out}(ii), and in view of lemma
\ref{lem_f-adic-pushout}, the quasi-affinoid scheme
$\sSpec\,\underline A$ represents the sought fibre
product (details left to the reader).
\end{example}

\begin{example}\label{ex_finite-ext-of-affinoids}
(i)\ \ 
Let $(A,A^+)$ be an affinoid ring and $B$ a finite
$A$-algebra. According to remark
\ref{rem_cantops-on-fin-algs}(iii), the canonical
topology $\cT^A_B$ on $B$ is f-adic, and the structure
map $\phi:A\to B$ is f-adic and restricts to a finite
map $A_0\to B_0$ of suitable subrings of definition.
Let $B^+$ be the integral closure of $\phi(A^+)$ in $B$.
We claim that $B^+$ is a ring of integral elements of
$B$. Indeed, $B^+\subset B^\circ$ by lemma
\ref{lem_f-adics}(iii.a) and remark
\ref{rem_someth-on-bdd-in-Z-lin}(iv). To check that $B^+$
is open in $B$, pick a finite system $b_1,\dots,b_n$ of
generators of the $A_0$-module $B_0$; for every
$i=1,\dots,n$ there exist $m_i\in\N$ and a polynomial
$P_i=X^{m_i}+\sum_{j=0}^{m_i-1}a_{ij}X^j\in A_0[X]$ with
$P_i(b_i)=0$. Let $I_0\subset A_0$ be an ideal of adic
definition; since $A^+$ is open, there exists $k\in\N$
such that $I^k_0a_{ij}\in A^+$ for $i=1,\dots,n$
and $j=0,\dots,m_i-1$. Thus, for every $c\in I^k_0$ and
$i=1,\dots,n$ we have
$$
0=c^{m_i}P_i(b_i)=Q_i(cb_i)
\qquad
\text{where $Q_i:=X^{m_i}+\sum_{j=0}^{m_i-1}c^{m_i-j}a_{ij}X^j\in A^+[X]$}.
$$
This shows that $I^k_0B_0\subset B^+$, whence the contention.
We call $(B,B^+)$ {\em the affinoid ring associated with the
finite $A$-algebra $B$}, and we denote it by
$$
B\otimes_A(A,A^+).
$$

(ii)\ \
Thus, we get an f-adic morphism $\phi:(A,A^+)\to B\otimes_A(A,A^+)$
of affinoid rings. We claim that the morphism of affinoid rings
$\phi$ enjoys the following universal property. For every morphism
$f:(A,A^+)\to(C,C^+)$ of affinoid rings, every ring homomorphism
$g:B\to C$ such that $g\circ\phi=f$ is a morphism
$g:(B,B^+)\to(C,C^+)$ of affinoid rings. Indeed, a simple
inspection shows that $g(B^+)\subset C^+$, and $g$ is continuous
by virtue of proposition \ref{prop_cantop-on-fg-mods}(ii).

(iii)\ \
As usual, the universal property determines the pair
$(B\otimes_A(A,A^+),\phi)$ up to unique isomorphism of
affinoid rings. Especially, if $\psi:B\to B'$ is any
finite ring homomorphism, then $\psi\circ\phi:A\to B'$
is also finite, and it follows easily that
$$
B'\otimes_B(B\otimes_A(A,A^+))=B'\otimes_A(A,A^+).
$$

(iv)\ \
Likewise, if $\underline A:=(A,A^+,U)$ is any quasi-affinoid
ring, and $\phi:A\to B$ is as in (i), let us set
$U_B:=\Spec\,B\times_{\Spec\,A}U$. According to example
\ref{ex_shouldbe-quasi-affinoid}(ii) we get a quasi-affinoid
ring
$$
B\otimes_A\underline A:=(B,B^+,U_B)
$$
(for the topology $\cT^A_B$ on $B$) and $\phi$ is an f-adic
morphism of quasi-affinoid rings
$\underline A\to B\otimes_A\underline A$ enjoying a
corresponding universal property that the reader may
spell out.

(v)\ \
In the same vein, let $\underline X:=(X,\cT_X,A^+_X)$ be a
quasi-affinoid scheme and $f:X'\to X$ a finite morphism of
schemes; set $A_X:=\cO_X(X)$, $Y:=\Spec\,A_X$ and let $i:X\to Y$
be the open immersion and $g:=i\circ f:X'\to Y$. We have
natural morphisms of quasi-coherent $\cO_Y$-algebras
$i^\flat:\cO_Y\to i_*\cO_X$ and $g^\flat:i_*\cO_X\to g_*\cO_{X'}$,
and we let $\cA$ be the integral closure of the image of
$i_*\cO_X$ in $g_*\cO_{X'}$. Notice that $i^\flat$ is an
isomorphism, and $\cA$ is a quasi-coherent $\cO_Y$-algebra;
hence $\cA$ is the filtered union of its finite quasi-coherent
$\cO_Y$-subalgebras.
Clearly $\cA_{|X}=f_*\cO_{X'}$ is a finite $\cO_X$-algebra;
since $X$ is quasi-compact, it follows that there exists a
finite quasi-coherent $\cO_Y$-subalgebra $\cB\subset\cA$ such
that $\cB_{|X}=f_*\cO_{X'}$, and $B:=\cB(Y)$ is a finite
$A_X$-algebra. We claim that the quasi-affinoid scheme
$$
X'\times_X\underline X:=\sSpec\,(B\otimes_{A_X}\sGamma(\underline X)).
$$
is independent, up to unique isomorphism, of the choice of
$\cB$ and enjoys the following universal property. For every
morphism of quasi-affinoid schemes
$h:\underline Z:=(Z,\cT_Z,A^+_Z)\to\underline X$, every morphism
$k:Z\to X'$ of schemes such that $f\circ k=h$ is a morphism of
quasi-affinoid schemes $\underline Z\to X'\times_X\underline X$.
Indeed notice that, by construction, the commutative diagram
of schemes
$$
\xymatrix{ X' \ar[r]^-j \ar[d]_f & \Spec\,B \ar[d] \\
X \ar[r]^-i & Y
}$$
is cartesian, hence $X'$ is the scheme underlying
$X'\times_X\underline X$. The morphisms $k$ and $j$ induce
homomorphisms of $A_X$-algebras $\cO_{X'}\to\cO_Z(Z)$ and
$B\to\cO_{X'}$; their composition is a homomorphism of
$A_X$-algebras $u:B\to\cO_Z(Z)$, and since $\cO_Z(Z)$ is
a topological $(A_X,\cT_X)$-algebra, (iv) implies that $u$
is a morphism of quasi-affinoid rings
$B\otimes_{A_X}\sGamma(\underline X)\to\sGamma(\underline Z)$.
By adjunction, $u$ induces a morphism of quasi-affinoid schemes
$\underline Z\to X'\times_X\underline X$, and a simple inspection
shows that this morphism is precisely given by $k$. Then the
universal property implies as usual the stated independence of
the choice of $\cB$. By the same token, we also deduce that
if $X''\to X'$ is another finite morphism of schemes, we have
a natural identification of quasi-affinoid schemes
$$
X''\times_{X'}(X'\times_X\underline X)\isom X''\times_X\underline X
$$
(details left to the reader). Lastly, notice that the
induced projection $X'\times_X\underline X\to\underline X$
is f-adic, by proposition \ref{prop_adjoint-to-sGamma}(ii).

(vi)\ \
In the situation of (i) (resp. of (iv), resp. of (v)),
suppose that $(A,A^+)$ (resp. $\underline A$, resp.
$\underline X$) is topologically henselian; combining
with proposition \ref{prop_quasi-affinoid}(iv) we deduce
that the same holds for $B\otimes_A(A,A^+)$ (resp. for
$B\otimes_A\underline A$, resp. for $X'\times_X\underline X$).
\end{example}

\begin{remark}\label{rem_loc-hens-cplete}
(i)\ \
Denote by $\mathsf{Afd.Ring}_\mathrm{comp}$ (resp.
$\mathsf{q.Afd.Ring}_\mathrm{comp}$) the full subcategory
of $\mathsf{Afd.Ring}$ (resp. of $\mathsf{q.Afd.Ring}$)
whose objects are the complete and separated affinoid
(resp. quasi-affinoid) rings. Then the inclusion functor
$$
\mathsf{Afd.Ring}_\mathrm{comp}\to\mathsf{Afd.Ring}
\qquad
\text{(resp.\ 
$\mathsf{q.Afd.Ring}_\mathrm{comp}\to\mathsf{q.Afd.Ring}$)}
$$
admits a left adjoint, that assigns to every affinoid ring
$(A,A^+)$ (resp. quasi-affinoid ring $(A,A^+,U)$) its
{\em completion}, which is the datum
$$
(A,A^+)^\wedge:=(A^\wedge,(A^+)^\wedge)
\qquad
\text{(resp. $(A,A^+,U)^\wedge:=(A^\wedge,(A^+)^\wedge,U^\wedge)$)}
$$
consisting of the separated completions of $A$ and $A^+$ (resp.
and of the open subset $U^\wedge:=U\times_{\Spec\,A}\Spec\,A^\wedge$).
Indeed, taking into account proposition
\ref{prop_complete-f-adic}(i,iii) and lemma
\ref{lem_compl-and-int.clos} we see that
$(A,A^+)^\wedge=A^\wedge\otimes_A(A,A^+)$ and
$(A,A^+,U)^\wedge=A^\wedge\otimes_A(A,A^+,U)$,
so the sought adjunctions follow from example
\ref{ex_shouldbe-quasi-affinoid}(i,ii).

(ii)\ \
Let $\mathsf{Afd.Ring}_\mathrm{hens}$ (resp.
$\mathsf{q.Afd.Ring}_\mathrm{hens}$) be the full subcategory
of $\mathsf{Afd.Ring}$ (resp. of $\mathsf{q.Afd.Ring}$)
whose objects are the topologically henselian affinoid
(resp. quasi-affinoid) rings. Then the inclusion functor
$$
\mathsf{Afd.Ring}_\mathrm{hens}\to\mathsf{Afd.Ring}
\qquad
\text{(resp.\ 
$\mathsf{q.Afd.Ring}_\mathrm{hens}\to\mathsf{q.Afd.Ring}$)}
$$
admits a left adjoint, that assigns to every affinoid
ring $(A,A^+)$ (resp. quasi-affinoid ring $(A,A^+,U)$)
its {\em topological henselization}, which is the datum
$$
(A,A^+)^\he:=(A^\he,A^{+\he})
\qquad
\text{(resp. $(A^\he,A^{+\he},U^\he)$)}
$$
consisting of the topological henselizations of $A$ and $A^+$
(resp. and of the open subset $U^\he:=U\times_{\Spec\,A}\Spec\,A^\he$ :
see definition \ref{def_henselize-f-adic}). Indeed, by
inspecting the construction of \eqref{subsec_top-henselization}
it is easily seen that $A^{+\he}$ identifies naturally
with an open subring of $A^\he$ contained in $(A^\he)^\circ$.
Moreover, $A^{+\he}$ is integrally closed in $A^\he$, by
corollary \ref{cor_pro-smooth}(ii). Hence $(A,A^+)^\he$ is
a topologically henselian affinoid ring, and theorem
\ref{th_hensel-f-adic}(ii) easily implies that every
morphism $(A,A^+)\to(B,B^+)$ of affinoid rings with
$(B,B^+)$ topologically henselian, factors uniquely
through a morphism $(A,A^+)^\he\to(B,B^+)$. By the same
token, we deduce that
$$
(A,A^+)^\he=A^\he\otimes_A(A,A^+)
$$
from which it follows that
$(A,A^+,U)^\he=A^\he\otimes_A(A,A^+,U)$, whence also
the second sought adjunction.

(iii)\ \
In the same vein, if $(A,A^+)$ (resp. $(A,A^+,U)$) is any
affinoid (resp. quasi-affinoid) ring, we may define its
{\em topological localization} by considering the datum
$(A,A^+)_\loc$ (resp. $(A,A^+,U)_\loc$)
consisting of the topological localizations $A_\loc$ and
$A^+_\loc$ of $A$ and respectively $A^+$, defined as in
\eqref{subsec_localize-f-adic} (resp. and of the open
subset $U\times_{\Spec\,A}\Spec\,A_\loc$). Then it is easily
seen that the rule $(A,A^+)\mapsto(A,A^+)_\loc$
(resp. $(A,A^+,U)\mapsto(A,A^+,U)_\loc$) yields
a left adjoint to the inclusion functor
$$
\mathsf{Afd.Ring}_\loc\to\mathsf{Afd.Ring}
\qquad
\text{(resp.\
$\mathsf{q.Afd.Ring}_\loc\to\mathsf{q.Afd.Ring}$)}
$$
from the full subcategory of $\mathsf{Afd.Ring}$
(resp. of $\mathsf{q.Afd.Ring}$) whose objects are the
topologically local affinoid (resp. quasi-affinoid) rings.
Notice also that $(A,A^+)_\loc=A_\loc\otimes_A(A,A^+)$,
and likewise for $(A,A^+,U)_\loc$.

(iv)\ \
Let $\underline X$ be any quasi-affinoid scheme; by definition,
$\underline X$ is topologically local (resp. topologically
henselian) if and only if the same holds for the quasi-affinoid
ring $\sGamma(\underline X)$.
On the other hand, let $\underline A:=(A,A^+,U)$ be any
quasi-affinoid ring; then $\underline A$ is topologically
local (resp. topologically henselian) if and only if the
same holds for the quasi-affinoid scheme $\sSpec\,\underline A$
(notation of proposition \ref{prop_adjoint-to-sGamma}).
Indeed, set $A_U:=\cO_U(U)$, and recall that the image $B$ of
$A$ in $A_U$ is open in the f-adic topology $\cT_U$ of $A_U$
(proposition \ref{prop_top-on-opens-fadic-case}(i)). In light
of proposition \ref{prop_quasi-affinoid}(i,iii), we are
then reduced to showing that $A$ is topologically local
(resp. topologically henselian) if and only if the same
holds for $B$, where the latter is endowed with the
f-adic topology induced by the inclusion map into $A_U$,
which is the same as the topology induced by $A$ via the
surjection $\pi:A\to B$. Let $U_A$ and $U_B$ be the analytic
loci of $\Spec\,A$ and respectively $\Spec\,B$; since $U$
contains $U_A$, it is easily seen that $\Spec\,\pi$ induces
a homeomorphism $U_B\isom U_A$ (details left to the reader).
On the other hand, notice that $B^{\circ\circ}$ is the
radical of the ideal $A^{\circ\circ}B$; then the assertion
follows from claim \ref{cl_Jacobson-or-Hensel} and
\cite[Rem.5.1.10(i)]{Ga-Ra}.

(v)\ \
For every quasi-affinoid scheme $\underline X$, set
$$
\underline X_\loc:=\sSpec(\sGamma(\underline X)_\loc)
\qquad
\underline X^\he:=\sSpec(\sGamma(\underline X)^\he).
$$
It follows easily from (iii),(iv), remark
\ref{rem_shouldbe-quasi-affinoid}(ii) and proposition
\ref{prop_adjoint-to-sGamma}(i) that $\underline X_\loc$
is a topologically local quasi-affinoid scheme endowed with
a natural morphism $j_{\underline X}:\underline X_\loc\to\underline X$
with the following property. Every morphism of quasi-affinoid
schemes $\underline Y\to\underline X$ with $\underline Y$
topologically local, factors uniquely through $j_{\underline X}$.
In other words, the rule
$\underline X\mapsto\underline X_\loc$ yields
a right adjoint
$$
\mathsf{q.Afd.Sch}\to\mathsf{q.Afd.Sch}_\loc
$$
to the inclusion functor of the full subcategory
$\mathsf{q.Afd.Sch}_\loc$ of topologically
local quasi-affinoid schemes. Likewise, $\underline X^\he$
is a topologically henselian quasi-affinoid scheme
endowed with a universal morphism
$h_{\underline X}:\underline X^\he\to\underline X$, so
the rule $\underline X\mapsto\underline X^\he$ yields
a right adjoint
$$
\mathsf{q.Afd.Sch}\to\mathsf{q.Afd.Sch}_\mathrm{hens}
$$
to the inclusion functor of the full subcategory
$\mathsf{q.Afd.Sch}_\mathrm{hens}$ of topologically
henselian quasi-affinoid schemes.

(vi)\ \
From (iv), we deduce that the adjunction $(\sSpec,\sGamma)$
of proposition \ref{prop_adjoint-to-sGamma}(i) restricts
to two adjoint pairs of functors
$$
\xymatrix{
\mathsf{q.Afd.Ring}_\loc \ar@<.5ex>[rr]^-{\sSpec} & &
\mathsf{q.Afd.Sch}^o_\loc \ar@<.5ex>[ll]^-{\sGamma} & 
\mathsf{q.Afd.Ring}_\mathrm{hens} \ar@<.5ex>[rr]^-{\sSpec} & &
\mathsf{q.Afd.Sch}^o_\mathrm{hens} \ar@<.5ex>[ll]^-{\sGamma}.
}$$
Taking into account (v) and remark \ref{rem_adjoint-transf}(i),
we  then obtain natural isomorphisms
$$
\sSpec(\underline A_\loc)\isom(\sSpec\,\underline A)_\loc
\qquad
\sSpec(\underline A_\mathrm{hens})\isom
(\sSpec\,\underline A)_\mathrm{hens}
$$
of quasi-affinoid schemes, for every quasi-affinoid ring
$\underline A$. On the other hand, we have as well natural
identifications
\set\begin{equation}\label{eq_wrong-order}
\sGamma(\underline X_\loc)\isom\sGamma(\underline X)_\loc
\qquad
\sGamma(\underline X^\he)\isom\sGamma(\underline X)^\he
\end{equation}
for every quasi-affinoid scheme $\underline X$. Indeed,
if $\sGamma(\underline X)=(A,A^+,U)$, we have
$\sGamma(\underline X)_\loc=(A_\loc,A^+_\loc,U_\loc)$, where
$U_\loc=U\times_{\Spec\,A}\Spec\,A_\loc$. Therefore
$\sGamma(\underline X_\loc)=(B,B^+,U_\loc)$, where
$B:=\cO_{\!U_\loc}(U_\loc)$, and with $B^+$ equal to the
integral closure of $A^+_\loc$ in $B$. Since the localization
map is flat, and since $A=\cO_{\!U}(U)$, the assertion for
$\sGamma(\underline X_\loc)$ then follows from corollary
\ref{cor_base-change-where}. The same argument applies
to $\sGamma(\underline X^\he)$.

(vii)\ \
On the other hand, if $\underline A:=(A,A^+,U)$ is a
complete and separated quasi-affinoid ring, the
quasi-affinoid scheme $\sSpec\,\underline A$ is not
necessarily complete and separated (essentially, this
is because the completion map $A\to A^\wedge$ is not
necessarily flat, in this generality). Thus, in order
to produce a completion functor for quasi-affinoid
schemes, we have to proceed more carefully, as explained
in the following proposition.
\end{remark}

\begin{proposition}\label{prop_complete-qaff-sch}
Let $\mathsf{q.Afd.Sch}_\mathrm{comp}$ be the full
subcategory of $\mathsf{q.Afd.Sch}$ whose objects
are the complete and separated quasi-affinoid schemes
$(X,\cT_X,A^+_X)$. Then the inclusion functor
$$
\mathsf{q.Afd.Sch}_\mathrm{comp}\to\mathsf{q.Afd.Sch}
$$
admits a right adjoint, called the {\em completion functor} :
$$
\mathsf{q.Afd.Sch}\to\mathsf{q.Afd.Sch}_\mathrm{comp}
\qquad
(X,\cT_X,A^+_X)\mapsto(X,\cT_X,A^+_X)^\wedge.
$$
\end{proposition}
\begin{proof} For any given quasi-affinoid scheme $(X,\cT_X,A_X^+)$,
set $A_X:=\cO_X(X)$, and denote by $(A_X^\wedge,\cT^\wedge_X)$ the
completion of $(A_X,\cT_X)$. For every ideal $I\subset A_X^\wedge$,
set as well
$$
V(I):=\Spec\,A_X^\wedge/I
\qquad\text{and}\qquad
X^\wedge_I:=X\times_{\Spec\,A_X}V(I).
$$
We let $\cF$ be the family of all topologically closed
ideals $I\subset A^\wedge_X$ such that the induced map
$\rho_I:A^\wedge_X/I\to\cO_{V(I)}(X^\wedge_I)$ is injective.
Set $J:=\bigcap_{I\in\cF}I$, so that the natural map
$i:A^\wedge_X/J\to\prod_{I\in\cF}A^\wedge_X/I$ is injective,
and notice that $(\prod_{I\in\cF}\rho_I)\circ i$ factors
naturally through $\rho_J$, therefore the latter is
injective, and thus $J\in\cF$. We endow $A':=A^\wedge_X/J$
with the quotient topology induced by the projection
$\pi:A^\wedge_X\to A'$; it is clear that $A'$ is a complete
and separated f-adic ring, and $\pi$ is f-adic (lemma
\ref{lem_f-adics}(iv)), so we may define
$$
(X,\cT_X,A^+_X)^\wedge:=
\sSpec\,(A'\otimes_{A^\wedge_X}(A_X,A^+_X,X)^\wedge).
$$
Since the restriction map $A'\to A'_X:=\cO_{V(J)}(X^\wedge_J)$
is injective and open, the topology of $A'_X$ is complete
and separated, so the datum $(X,\cT_X,A^+_X)^\wedge$ is indeed
an object of $\mathsf{q.Afd.Sch}_\mathrm{comp}$.

On the other hand, let $(Y,\cT_Y,A^+_Y)$ be any complete
and separated quasi-affinoid scheme, and
$f:(Y,\cT_Y,A^+_Y)\to(X,\cT_X,A^+_X)$ any morphism of
quasi-affinoid schemes; by adjunction, $\sGamma(f)$
factors uniquely through a morphism of complete
quasi-affinoid rings
$$
\sGamma(f)^\wedge:(A_X,A^+_X,X)^\wedge\to\sGamma(Y,\cT_Y,A^+_Y)
$$
and we let $K$ be the kernel of the underlying
continuous ring homomorphism $h:A^\wedge_X\to A_Y$.
Clearly $K$ is a closed ideal of $A_X$, and $\Spec\,h$
restricts to a morphism of schemes $g:Y\to V(K)$; then,
the pair $(f,g)$ determines a unique morphism of
$X$-schemes $Y\to X^\wedge_K$. Therefore, the resulting
map $\bar h:A^\wedge_X/K\to A_Y$ is injective by
construction, and it factors through $\rho_K$;
so $\rho_K$ is injective as well, and thus $K\in\cF$.
Hence, $h$ factors uniquely through a continuous
ring homomorphism $A'\to A_Y$, so in turns $\sGamma(f)^\wedge$
factors uniquely through a morphism of quasi-affinoid
rings
$$
\phi:A'\otimes_{A_X}(A_X,A^+_X,X)^\wedge\to\sGamma(Y,\cT_Y,A^+_Y).
$$
and $\sSpec(\phi):(Y,\cT_Y,A^+_Y)\to(X,\cT_X,A^+_X)^\wedge$
is the unique morphism of complete quasi-affinoid schemes
that lifts $f$, whence the proposition.
\end{proof}

\begin{remark}\label{rem_depth-and-completion}
(i)\ \
Let $\underline X:=(X,\cT_X,A^+_X)$ be any quasi-affinoid scheme,
set $(A_X,A^+_X,X):=\sGamma(\underline X)$, and let $A_X^\wedge$
be the separated completion of $A_X$. Let also $I\subset A_X$
be any finitely generated ideal such that
$\Spec\,A_X/I=\Spec\,A_X\setminus X$. By inspecting the
proof of proposition \ref{prop_complete-qaff-sch} we obtain
a natural morphism of quasi-affinoid schemes
$$
\underline X^\wedge\to\sSpec(\sGamma(\underline X)^\wedge)
$$
which is an isomorphism if and only if $\depth_IA_X^\wedge>0$
(notation of \eqref{subsec_depth_A}).

(ii)\ \
The morphism of (i) induces by adjunction a morphism of
quasi-affinoid rings
$$
\sGamma(\underline X)^\wedge\to\sGamma(\underline X^\wedge)
$$
which is an isomorphism if and only if $\depth_IA_X^\wedge>1$.
\end{remark}

\begin{definition}\label{def_adic-spectrum}
Let $\underline A:=(A,A^+,U)$ be any quasi-affinoid ring,
and $\underline X:=(X,\cT_X,A^+_X)$ any quasi-affinoid scheme.

(i)\ \
The {\em adic spectrum} of the affinoid ring $(A,A^+)$
is the set
$$
\Spa(A,A^+):=
\{v\in\Cont(A)~|~\text{$v(a)\leq 1$ for every $a\in A^+$}\}
$$
(notation of definition \ref{def_continuous-vals}).
Hence, $\Spa(A,A^+)$ is a pro-constructible subset of
$\Cont(A)$, and we endow $\Spa(A,A^+)$ with the topology
induced by $\Cont(A)$. We also let
$$
\Spa(A,A^+)_\mathrm{a}:=\Spa(A,A^+)\cap\Cont(A)_\mathrm{a}
\qquad
\Spa(A,A^+)_\mathrm{na}:=\Spa(A,A^+)\cap\Cont(A)_\mathrm{na}
$$
(notation of remark \ref{rem_Cont-A}(vi)).

(ii)\ \
The {\em adic spectrum} of $\underline A$ is the set
$$
\Spa\,\underline A:=\Spa(A,A^+)\cap\Spv\,U.
$$
We endow this subset with the topology induced by
$\Spa(A,A^+)$.

(iii)\ \
The {\em adic spectrum} of $\underline X$ is the topological
space
$$
\Spa\,\underline X:=\Spa\,\sGamma(\underline X)
$$
and we set as well
$$
\Spa\,(\underline X)_\mathrm{a}:=
\Spa\,\sGamma(\underline X)_\mathrm{a}
\qquad\text{and}\qquad
\Spa\,(\underline X)_\mathrm{na}:=
\Spa\,\sGamma(\underline X)_\mathrm{na}.
$$

(iv)\ \
A {\em rational subset} of $\Spa(A,A^+,U)$ is a subset
of the form $R\cap\Spa(A,A^+,U)$, where $R$ is a rational
subset of $\Cont(A)$ (see \eqref{subsec_rational-Cont}).
Likewise, a {\em rational subset} of $\Spa\,\underline X$
is a subset of the form $R\cap\Spa\,\underline X$, where
$R$ is a rational subset of $\Cont(\cO_{\!X}(X),\cT_X)$.
\end{definition}

\begin{remark}\label{rem_about-Spa}
(i)\ \ 
Let $f:(A,A^+)\to(B,B^+)$ be any morphism of affinoid rings.
A simple inspection of definition \ref{def_adic-spectrum}(i)
shows that the map $\Cont(f)$ restricts to a continuous map
$$
\Spa\,f:\Spa(B,B^+)\to\Spa(A,A^+).
$$

(ii)\ \
More generally, a morphism $g:(A,A^+,U)\to(B,B^+,V)$ of
quasi-affinoid rings induces a continuous map
$$
\Spa\,g:\Spa(B,B^+,V)\to\Spa(A,A^+,U)
$$
that agrees with the restriction of the map
$\Spa(B,B^+)\to\Spa(A,A^+)$ attached to $g$ as in (i).

(iii)\ \
Likewise, any morphism $\phi:(X,\cT_X,A^+_X)\to(Y,\cT_Y,A^+_Y)$
of quasi-affinoid schemes induces a continuous map
$$
\Spa\,\phi:=\Spa\,\sGamma(\phi):
\Spa(X,\cT_X,A^+_X)\to\Spa(Y,\cT_Y,A^+_Y).
$$

(iv)\ \
In light of theorem \ref{th_Cont-spectral}(i), it is easily
seen that $v(x)\leq 1$ for every $x\in A^\circ$ and every
$v\in\Cont(A)$ of rank $\leq 1$. Especially, the subset
of $\Spa(A,A^+,U)$ consisting of rank one valuations is
independent of the ring of integral elements $A^+$.
\end{remark}

\begin{example}\label{ex_affinoids}
For any f-adic ring $A$, let $A'\subset A$ be the
smallest subring of integral elements (see remark
\ref{rem_shouldbe-quasi-affinoid}(iv)). It is easily
seen that $\Spa(A,A')=\Cont(A)$. Also, every continuous
map $f:A\to B$ of f-adic rings yields a morphism of
affinoid rings $f:(A,A')\to(B,B')$, and clearly
$\Spa\,\phi=\Cont(f)$.

(ii)\ \
On the other hand, every f-adic map $f:A\to B$ of
f-adic rings induces a continuous map
$\Spa\,f:\Spa(B,B^\circ)\to\Spa(A,A^\circ)$ (lemma
\ref{lem_f-adics}(iii.a)).
\end{example}

\begin{lemma}\label{lem_invariance-by-loc-hens}
Let $\underline A:=(A,A^+,U)$ be any quasi-affinoid ring.
\begin{enumerate}
\item
The natural morphisms
$\underline A\to\underline A{}_\loc\to\underline A^\he\to
\underline A^\wedge$
(notation of remark {\em\ref{rem_loc-hens-cplete}(ii,iii)})
induce homeomorphisms
$$
\Spa\,\underline A^\wedge\isom\Spa\,\underline A^\he\isom
\Spa\,\underline A{}_\loc\isom\Spa\,\underline A.
$$
\item
The unit of the adjunction $(\sSpec,\sGamma)$ of proposition
{\em\ref{prop_adjoint-to-sGamma}(i)} induces a homeomorphism
$$
\Spa(\sSpec\,\underline A)\isom\Spa\,\underline A.
$$
\item
Let $\underline X:=(X,\cT_X,A^+_X)$ be any quasi-affinoid
scheme, and $\underline X^\wedge$ the completion of
$\underline X$ (see proposition {\em\ref{prop_complete-qaff-sch}}).
The completion map $\underline X\to\underline X^\wedge$ induces
a homeomorphism
$$
\Spa\,\underline X^\wedge\isom\Spa\,\underline X.
$$
\item
Let $\phi:\Spa\,\underline A^\wedge\isom\Spa\,\underline A$
be the homeomorphism of\/ {\em (i)}, and
$R\subset\Spa\,\underline A^\wedge$ any rational subset.
Then $\phi(R)$ is a rational subset of\/ $\Spa\,\underline A$.
\end{enumerate}
\end{lemma}
\begin{proof}(i): Let $A_\loc$ (resp. $A^\he$) be the
topological localization (resp. henselization) of $A$;
from corollary \ref{cor_justify}(i) and proposition
\ref{prop_Cont-of-complete}(iv), we deduce that the
natural maps $A\to A_\loc\to A^\he\to A^\wedge$
induce homeomorphisms
$$
\Cont(A^\wedge)\isom\Cont(A^\he)\isom\Cont(A_\loc)
\isom\Cont(A).
$$
Moreover, it was already remarked that the localization
map induces a natural identification
$\underline A{}_\loc\isom A_\loc\otimes_A(A,A^+)$,
and likewise for the topological henselization and
the completion of $\underline A$; the assertion follows
easily (details left to the reader).

(ii): Set $(U,\cT_U,A^+_U):=\sSpec\,\underline A$ and
$A_U:=\cO_U(U)$; the assertion boils down to the following:

\begin{claim}
The natural morphism of quasi-affinoid rings
$\rho:(A,A^+,U)\to(A_U,A^+_U,U)$ induces a homeomorphism
$\Spa\,\rho:\Spa(A_U,A^+_U,U)\isom\Spa(A,A^+,U)$.
\end{claim}
\begin{pfclaim} Let us first check that $\rho$ induces
a homeomorphism
\set\begin{equation}\label{eq_this-map}
\Cont(A_U)\cap\Spv\,U\isom\Cont(A)\cap\Spv\,U.
\end{equation}
Indeed, since the topology of both of these spaces are
induced from the inclusion into $\Spv\,U$, it suffices
to check that \eqref{eq_this-map} is surjective. Hence,
let $v:A\to\Gamma_{\!v\circ}$ be any continuous valuation
with support given by a prime ideal $\fp\subset A$, and
suppose that $\fp\in U$; then there exists a unique
prime ideal $\fp'\subset A_U$ such that $\rho^{-1}\fp'=\fp$.
Set $B:=A/\fp$ and $C:=A_U/\fp'$; the map $v$ factors
through a continuous valuation $\bar v:B\to\Gamma_{\!v\circ}$
(for the quotient topology on $B$ induced by the projection
$A\to B$) and $C$ is naturally identified with a subring
of $\Frac\,B$, so $\bar v$ extends uniquely to a
valuation $w:C\to\Gamma_{\!v\circ}$. It remains only
to show that $w$ is continuous, for the quotient topology
on $C$ induced by the projection $A_U\to C$. By remark
\ref{rem_Cont-A}(ii) it suffices to prove the continuity of
$w$ at the point $0\in C$; however, $\rho$ is an open map
(proposition \ref{prop_top-on-opens-fadic-case}(i)), so the
same holds for the induced injective map $B\to C$, and the
assertion follows.

Since $A^+_U$ is the integral closure of the image of $A^+$
in $A_U$, the claim is an immediate consequence of the
foregoing : details left to the reader.
\end{pfclaim}

(iii): Let $A_X:=\cO_{\!X}(X)$, and set
$(A_X^\wedge,A_X^{+\wedge},X^\wedge):=(A_X,A^+_X,X)^\wedge$.
Define the family $\cF$ of ideals of $A^\wedge_X$ as in
the proof of proposition \ref{prop_complete-qaff-sch},
and denote by $J$ the minimal element of $\cF$; also,
endow $A':=A^\wedge_X/J$ with the quotient topology
induced via the projection $A^\wedge_X\to A'$. Then
$\underline X^\wedge=
\sSpec(A'\otimes_{A^\wedge_X}(A_X^\wedge,A_X^{+\wedge},X^\wedge))$.
Taking into account (i) and (ii), we are then reduced
to checking that the natural projection
$\pi:(A_X^\wedge,A_X^{+\wedge},X^\wedge)\to
A'\otimes_{A^\wedge_X}(A_X^\wedge,A_X^{+\wedge},X^\wedge)$ induces
a homeomorphism
$$
\Spa\,\pi:
\Spa(A'\otimes_{A^\wedge_X}(A_X^\wedge,A_X^{+\wedge},X^\wedge))\isom
\Spa(A_X^\wedge,A_X^{+\wedge},X^\wedge).
$$
However, $\Spa\,\pi$ is a closed immersion, so it suffices
to show that it is a bijection. To this aim, let
$v:A^\wedge_X\to\Gamma_{\!v\circ}$ be any continuous
valuation whose support $\fp$ lies in $X^\wedge$; we come
down to checking that $\fp\in\cF$. However, $\fp$ is
obviously a topologically closed ideal of $A^\wedge_X$.
Next, set $B:=A^\wedge_X/\fp$ and
$X^\wedge_\fp:=X\times_{\Spec\,A_X}\Spec\,B$; since $B$ is a
domain and $X^\wedge_\fp\neq\emptyset$, the restriction map
$B\to\cO_{\Spec\,B}(X^\wedge_\fp)$ is injective, whence the
contention.

(iv): Say that
$R=R_{A^\wedge}\bigl(\frac{a_1}{a_0},\dots,\frac{a_n}{a_0}\bigr)
\cap\Spa\,\underline A^\wedge$ for a sequence of elements
$a_0,\dots,a_n\in A^\wedge$ that generates an open ideal $J$.
Let also $B$ be a ring of definition of $A$, and
$I\subset B$ an ideal of adic definition. To begin with,
we remark :

\begin{claim}\label{cl_efficient}
(i)\ \
Let $T$ be a quasi-compact subset of $\Spa\,\underline A$,
and $t$ an element of $A^{\circ\circ}$ such that $v(t)\neq 0$
for every $v\in T$. Then there exists $k\in\N$ such that
$v(x)<v(t)$ for every $v\in T$ and every $x\in I^k$.

(ii)\ \
There exists $q\in\N$ such that every sequence
$b_0,\dots,b_n$ with $b_i\in a_i+I^qB^\wedge$ for
$i=0,\dots,n$, generates an open ideal of $A^\wedge$, and
$R':=R_{A^\wedge}\bigl(\frac{b_1}{b_0},\dots,\frac{b_n}{b_0}\bigr)
\cap\Spa\,\underline A^\wedge=R$.
\end{claim}
\begin{pfclaim}(i): Let $x_1,\dots,x_m$ be a finite
system of generators of $I$. For every $r\in\N$ set
$$
T_r:=\{v\in\Spa\,\underline A~|~v(x_i^r)\leq v(t)\ 
\text{for every $i=1,\dots,m$}\}.
$$
Then $T_r$ is open in $\Spa\,\underline A$ for every
$r\in\N$, and $T\subset\bigcup_{r\in\N}T_r$. Thus,
$T\subset T_r$ for some $r\in\N$, and the assertion
holds with $k:=r+1$.

(ii): Pick elements $c_1,\dots,c_m\in J\cap B^\wedge$
that generate an open ideal $J_0$ of $B^\wedge$, and $r\in\N$
such that $I^r\subset J_0$. Then, every sequence $d_1,\dots,d_m$
of elements of $B^\wedge$ with $d_i\in c_i+I^rB^\wedge$ for
$i=1,\dots,m$, also generates $J_0$ (\cite[Th.8.4]{Mat}).
Next, for every $i=1,\dots,m$ pick $x_{i0},\dots,x_{in}\in A^\wedge$
such that $c_i=\sum_{j=0}^na_jx_{ij}$, and let $p\in\N$ be large
enough so that $I^px_{ij}\subset I^rB^\wedge$ for every $i=1,\dots,m$
and every $j=0,\dots,n$; it follows already that every sequence
$b_0,\dots,b_n$ of elements of $A^\wedge$ with $b_i\in a_i+I^pB^\wedge$
for $i=0,\dots,n$ generates an open ideal of $A^\wedge$.

Next, we apply part (i) of the claim to the rational subsets
$R_i:=R_{A^\wedge}\bigl(\frac{a_1}{a_i},\dots,\frac{a_n}{a_i}\bigr)
\cap\Spa\,\underline A^\wedge$, to find an integer $q\geq p$
such that $v(x)<v(a_i)$ for every $x\in I^qB^\wedge$, every
$i=0,\dots,n$ and every $v\in R_i$. We claim that this $q$
will do. Indeed, let $v\in\Spa\,\underline A^\wedge$, and say
that $v\in R=R_0$; since $b_i-a_i\in I^qB^\wedge$, we have
$v(b_i-a_i)<v(a_0)$ for $i=0,\dots,n$. This implies that
$$
v(b_i)=v(a_i+(b_i-a_i))\leq\max(v(a_i),v(b_i-a_i))\leq
v(a_0)=v(a_0+(b_0-a_0))=v(b_0)
$$
for every $i=0,\dots,n$, which shows that $R\subset R'$.
Next, suppose $v\notin R$, and consider first the case
where $v(a_i)=0$ for every $i=0,\dots,n$. Then $v$ is
non-analytic, so that $v(b_0-a_0)=0$, whence $v(b_0)=0$,
and finally $v\notin R'$. Lastly, suppose that $v(a_i)\neq 0$
for some $i\leq n$, and choose $j\leq n$ so that $v\in R_j$.
We must have $v(a_0)<v(a_j)$, since $v\notin R$; also, by
construction $v(b_i-a_i)<v(a_j)$ for every $i=0,\dots,n$.
Then
$$
v(b_0)=v(a_0+(b_0-a_0))\leq\max(v(a_0),v(b_0-a_0))<
v(a_j)=v(a_j+(b_j-a_j))=v(b_j)
$$
which yields again $v\notin R'$, and the claim follows.
\end{pfclaim}

Now, by claim \ref{cl_efficient}(ii) there exist $b_0,\dots,b_n$
in $A$ such that
$R=R_{A^\wedge}\bigl(\frac{b_1}{b_0},\dots,\frac{b_n}{b_0}\bigr)\cap
\Spa\,\underline A^\wedge$. Also, by claim \ref{cl_efficient}(i)
there exists $k\in\N$ such that $v(x)<v(b_0)$ for every $v\in R$
and every $x\in I^k$.
Pick a finite system $b_{n+1},\dots,b_p$ of generators of $I^k$;
it follows that
$\phi(R)=R_A\bigl(\frac{b_1}{b_0},\dots,\frac{b_p}{b_0}\bigr)$
whence the contention.
\end{proof}

\begin{remark} By combining lemma
\ref{lem_invariance-by-loc-hens}(iv) and proposition
\ref{prop_Cont-of-complete}(i) one may obtain another,
more constructive, proof of lemma
\ref{lem_invariance-by-loc-hens}(i) (and of proposition
\ref{prop_Cont-of-complete}(iv)). This is the route
followed in \cite{Hu2}. 
\end{remark}

\begin{proposition}\label{prop_f-adic-Spectra}
Let $(A,A^+,U)$ be any quasi-affinoid ring. The following
holds :
\begin{enumerate}
\item
$\Spa(A,A^+,U)$ is a constructible open subset of\/
$\Spa(A,A^+)$, and both these topological spaces
are spectral.
\item
The rational subsets are a basis of quasi-compact open
subsets of\/ $\Spa(A,A^+)$ that is closed under finite
intersections.
\item
For any morphism $f:(A,A^+,U)\to(B,B^+,V)$ of
quasi-affinoid rings, we have :
\begin{enumerate}
\item
If $f$ is f-adic, $\Spa\,f$ is a spectral map. More
precisely, if $R$ is a rational subset of\/ $\Spa(A,A^+,U)$,
then $(\Spa\,f)^{-1}(R)$ is a rational subset of\/
$\Spa(B,B^+,V)$.
\item
Suppose that $B$ is topologically local. Then $f$
is f-adic if and only if\/ $\Spa\,f$ restricts to a map
$\Spa(B,B^+)_\mathrm{a}\to\Spa(A,A^+)_\mathrm{a}$.
\end{enumerate}
\item
For any quasi-affinoid scheme $(X,\cT_X,A^+_X)$, the
topological space $\Spa(X,\cT_X,A^+_X)$ is spectral.
\item
For any f-adic morphism $\phi:(X,\cT_X,A^+_X)\to(Y,\cT_Y,A^+_Y)$
of quasi-affinoid schemes, the induced map $\Spa\,\phi$
is spectral.
\end{enumerate}
\end{proposition}
\begin{proof} Assertion (i) for $\Spa(A,A^+)$ follows
from theorem \ref{th_Cont-spectral}(ii) and corollary
\ref{cor_procon-is-spec}. According to
\eqref{subsec_rational-Cont}, the rational subsets of
$\Cont(A)$ form a basis of quasi-compact open subsets
for the latter space; invoking again corollary
\ref{cor_procon-is-spec}, we obtain (ii) as well.

Next, let $I\subset A$ be any finitely generated ideal
such that $\Spec\,A/I=\Spec\,A\setminus U$, and pick a
finite system of generators $a_1,\dots,a_n$ for $I$. Then
$$
\Spa(A,A^+,U)=\Spa(A,A^+)\cap
\bigcup_{i=1}^nR_A\Bigl(\frac{a_1}{a_i},\dots,\frac{a_n}{a_i}\Bigr).
$$
Since $I$ is open in $A$, this identity presents $\Spa(A,A^+,U)$
as a finite union of rational subsets of $\Spa(A,A^+)$, and
the latter are open and quasi-compact, by (ii), so the proof
of (i) is complete.

(iii.a) follows immediately from corollary
\ref{cor_Cont-spectral}(i). The proofs of (iv) and (v)
are similar, and the details shall be left to the reader.

(iii.b): Set $\phi:=\Spec\,f$; in light of lemma
\ref{lem_deja-vu}(iv) and corollary \ref{cor_characterize-f-adic},
we need only check that $\phi^{-1}X^{\circ\circ}_A=X^{\circ\circ}_B$,
if $\Spa\,f$ restricts to a map
$\Spa(B,B^+)_\mathrm{a}\to\Spa(A,A^+)_\mathrm{a}$.

Thus, suppose that there exists
$\fp\in X_B\setminus X^{\circ\circ}_B$ with
$\phi(\fp)\in X_A^{\circ\circ}$. Then $\fp$ is a non-open
ideal of $B$, so there exists a rank one continuous
analytic valuation $v$ of $B$ whose support contains
$\fp$ (lemma \ref{lem_criterion-opennes}); clearly the
support of $\Cont(f)(v)$ contains $\phi(\fp)$, and it
is therefore open, {\em i.e.}
$\Cont(f)(v)\in\Cont(A)_\mathrm{na}$. To conclude
it suffices to remark :

\begin{claim}\label{cl_rank-one-vals}
Let $R$ be any f-adic ring, $v\in\Cont(R)_\mathrm{a}$
any element such that $\Gamma_{\!v}$ has rank one. Then
$v\in\Spa(R,R^\circ)$.
\end{claim}
\begin{pfclaim}[] Since $v$ is analytic, there exists
$a\in R$ such that $0<v(a)<1$. Let $b\in R^\circ$ be
any element, and suppose that $v(b)>1$; then there
exists $n\in\N$ large enough, such that $v(b^na)>1$,
which is absurd, since $b^na\in R^{\circ\circ}$ (remark
\ref{rem_someth-on-bdd-in-Z-lin}(iv)).
\end{pfclaim}
\end{proof}

\begin{proposition}\label{prop_technical-prop}
Let $(A,A^+,U)$ be any quasi-affinoid ring, $I\subset A^+$
an ideal, $g_\bullet:=(g_i~|~i\in\Sigma)$ and
$h_\bullet:=(h_j~|~j\in\Sigma')$ two systems of elements of
$A$ such that :
\begin{enumerate}
\alphaenu
\item
The system $(g_\bullet,h_\bullet)$ generates an open
ideal of $A$.
\item
$\Spec\,A/IA$ is the set of open prime ideals of $A$.
\item
$U$ is dense in $\Spec\,A$.
\end{enumerate}
Then for every $f\in A$ the following conditions are
equivalent :
\begin{enumerate}
\addenu\addenu\addenu
\alphaenu
\item
For every $v\in\Spa(A,A^+,U)$ with $v(f)\neq 0$ there
exists either $i\in\Sigma$ such that $v(f)\leq v(g_i)$
or else $j\in\Sigma'$ such that $v(f)<v(h_j)$.
\item
There exists a polynomial $P\in\cP_I$ such that
$P(f,g_\bullet,h_\bullet)=0$ (notation of \eqref{subsec_technical}).
\end{enumerate}
\end{proposition}
\begin{proof} From theorem \ref{th_Cont-spectral}(i) and
lemmata \ref{lem_invariance-by-loc-hens}(i),
\ref{lem_technical-lemma} we already see that
(e)$\Rightarrow$(d). 

(d)$\Rightarrow$(e): Let $(A',A'^+,U')$ be the topological
localization of $(A,A^+,U)$, and set $I':=IA'^+$. Notice
that $U'$ is dense in $\Spec\,A'$ (proposition
\ref{prop_closed-under-spec}(iii)), and $\Spec\,A'/I'A'$ is
the set of open prime ideals of $A'$, since the localization
map $A\to A'$ is f-adic (lemma \ref{lem_deja-vu}(iv)). We
may then consider also the set
$\cP_{I'}\subset A'^+[Z,X_\bullet,Y_\bullet]$ defined as in
\eqref{subsec_technical}, and we remark

\begin{claim}\label{cl_if-e-fails}
Condition (e) is equivalent to :
\begin{itemize}
\item[(f)]
There exists a polynomial $Q\in\cP_{I'}$ such that
$Q(f,g_\bullet,h_\bullet)=0$ in $A'^+$.
\end{itemize}
\end{claim}
\begin{pfclaim} Obviously (e)$\Rightarrow$(f). Conversely,
let $A_0\subset A^+$ be a subring of definition, and $J$ any
ideal of adic definition of $A_0$, so that $A'^+=(1+J)^{-1}A^+$
and $I'=(1+J)^{-1}I$. If (f) holds, we deduce that there
exists a polynomial $G\in A^+[Z,X_\bullet,Y_\bullet]$ of the form
$$
G(Z,X_\bullet,Y_\bullet)=
(1-a)\cdot Z^n+\sum_{k=1}^nZ^{n-k}P_k(X_\bullet,Y_\bullet)
$$
for some $a\in J$ and some $P_1,\dots,P_n\in A^+[X_\bullet,Y_\bullet]$
fulfilling the conditions of \eqref{subsec_technical}, and
such that $G(f,g_\bullet,h_\bullet)=0$ in $A$. Since $A$ is
f-adic, condition (a) implies that the system
$(g_i^n,h_j^n~|~i\in\Sigma,\ j\in\Sigma')$ generates
an open ideal of $A$ (details left to the reader);
then we may find $N\in\N$ large enough, finite subsets
$\Sigma_0\subset\Sigma'$, $\Sigma'_0\subset\Sigma'$, and
systems $(b_i~|~i\in\Sigma_0)$, $(c_j~|~j\in\Sigma'_0)$
of elements of $A$ such that
$$
a^Nf^n=\sum_{i\in\Sigma_0}b_ig_i^n+\sum_{j\in\Sigma'_0}c_jh^n_j.
$$
In light of (b), we may next find $M\in\N$ sufficiently
large, so that $a^Mb_i\in A^+$ and $a^Mc_j\in I$ for every
$i\in\Sigma_0$ and $j\in\Sigma'_0$. Summing up, we see that
$a^{M+N}f^n=H(g_\bullet,h_\bullet)$ in $A$ for some homogeneous
polynomial $H(X_\bullet,Y_\bullet)\in A^+[X_\bullet,Y_\bullet]$
of degree $n$ and with $H(0,Y_\bullet)\in I[Y_\bullet]$.
Then set $u:=1+a+\cdots+a^{M+N-1}$, so that
$u\cdot(1-a)=1-a^{M+N}$; it follows that the polynomial
$$
P:=u\cdot G+a^{M+N}Z^n-H
$$
fulfills condition (e).
\end{pfclaim}

Now, if (e) fails for some $f\in A$, claim \ref{cl_if-e-fails}
and lemma \ref{lem_technical-lemma} show that there
exists $w\in\Spv\,A'$ whose support is a minimal prime
ideal of $A'$, and with $w(f)\neq 0$ and
$$
\begin{aligned}
w(a)\leq &\, 1 \qquad & w(b)< &\, 1 \qquad & &
\text{for every $a\in A'^+$ and $b\in A'^{\circ\circ}$} \\
w(g_i)< &\, w(f) \qquad & w(h_j)\leq &\, w(f) \qquad & &
\text{for every $i\in\Sigma$ and $j\in\Sigma'$}.
\end{aligned}
$$
Since $U'$ is dense in $\Spec\,A'$, we have $w\in\Spv\,U'$
(proposition \ref{prop_closed-under-spec}(iii)); if $w$ is
continuous, it follows that $w\in\Spa(A',A'^+,U')=\Spa(A,A^+,U)$
(lemma \ref{lem_invariance-by-loc-hens}(i)), and (d) fails for
$w$. In case $w$ is not continuous, set $L:=IA$ and
$v:=w^{c\Gamma_v(L)}$ (notation of definition \ref{def_cGamma-I}).

\begin{claim}\label{cl_w-is-ok}
(i)\ \ 
$v(L)\neq\{0\}$.
\begin{enumerate}
\addenu
\item
$v(f)\neq 0$.
\end{enumerate}
\end{claim}
\begin{pfclaim}(i): Since $v$ is an $L$-admissible
specialization of $w$, it suffices to show that
$w(L)\neq\{0\}$. But if the latter fails, then the
support of $w$ is an open prime ideal, so $w$ is
continuous (remark \ref{rem_Cont-A}(vi)),
contradicting our assumption.

(ii): Suppose that $v(f)=0$; then $v(g_i)=v(h_j)=0$
as well, for every $i\in\Sigma$ and $j\in\Sigma'$.
Hence the support of $v$ is an open ideal of $A$,
and especially, it contains $L$, contradicting (i).
\end{pfclaim}

Claim \ref{cl_w-is-ok}(ii), theorem \ref{th_Cont-spectral}(i)
and lemma \ref{lem_invariance-by-loc-hens}(i) imply that
$v$ is an analytic point of $\Spa(A,A^+)$ and condition
(d) fails for $v$; lastly, notice that
$\Cont(A)_\mathrm{a}\subset\Spv\,U$, so $v\in\Spa(A,A^+,U)$
whence the proposition.
\end{proof}

\begin{corollary}\label{cor_corcor}
Let $\underline A:=(A,A^+,U)$ be any quasi-affinoid ring such
that $U$ is dense in $\Spec\,A$, and $f\in A$ any element.
We have :
\begin{enumerate}
\item
$f\in A^+$ if and only if $v(f)\leq 1$ for every
$v\in\Spa\,\underline A$.
\item
$f\in A^{\circ\circ}$ if and only if $v(f)<1$ for every
$v\in\Spa\,\underline A$.
\item
Let also $J\subset A^+$ be any finitely generated ideal
whose radical is $A^{\circ\circ}$. Then the following
conditions are equivalent :
\begin{enumerate}
\item
$|f|^*_J=0$ (notation of example {\em\ref{ex_Samuel}}).
\item
$v(f)=0$ for every $v\in\Spa\,\underline A$.
\end{enumerate}
\end{enumerate}
\end{corollary}
\begin{proof}(i): We may assume that $v(f)\leq 1$ for every
$v\in\Spa\,\underline A$, and we show that $f\in A^+$.
To this aim, we apply proposition \ref{prop_technical-prop}
with $g_\bullet=\{1\}$, $h_\bullet$ the empty subset and
$I:=A^{\circ\circ}$; we deduce that there exists a monic
polynomial $P(Z)\in A^+[Z]$ with $P(f)=0$, {\em i.e.}
$f$ is integral over $A^+$, whence the contention, since
$A^+$ is integrally closed in $A$.

(ii): Again, we may assume that $v(f)<1$ for every
$v\in\Spa\,\underline A$, and we check that $f\in A^{\circ\circ}$.
To this aim, we apply proposition \ref{prop_technical-prop}
with $h_\bullet=\{1\}$, $g_\bullet$ the empty subset and
$I:=A^{\circ\circ}$; we deduce that there exist $n\in\N$
and $a_1,\dots,a_n\in A^{\circ\circ}$ with
$f^n+f^{n-1}a_1+\cdots+a_n=0$. By (i), we know already that
$f\in A^+$, so that $f^{n-1}a_1+\cdots+a_n\in A^{\circ\circ}$,
and combining with remark
\ref{rem_someth-on-bdd-in-Z-lin}(iv), the claim follows.

(iii.a)$\Rightarrow$(iii.b): Let $a_1,\dots,a_r$ be any
finite system of generators of $J$, and set
$$
\gamma_v:=\max(v(a_1),\dots,v(a_r))\in\Gamma_{\!v\circ}
\qquad
\text{for every $v\in\Spa(A,A^+,U)$}.
$$
The assumption implies that for every $n\in\N$ there
exists $k\in\N$ such that $f^k\in J^{kn}$, hence
$v(f^k)\leq\gamma_v^{kn}$, and therefore
$v(f)\leq\gamma_v^n$ for every $n\in\N$ and every
$v\in\Spa\,\underline A$ (details left to the reader).
However $\gamma_v$ is final in $\Gamma_{\!v}$ for
every such $v$ (lemma \ref{lem_Cont-A}(i)), whence the
assertion.

(iii.b)$\Rightarrow$(iii.a): For every $n\in\N$,
we apply the criterion of proposition
\ref{prop_technical-prop} with $g_\bullet$ equal to the
empty set, $h^{(n)}_\bullet:=\{1\}$, and $I:=J^n$. We
deduce that for every $n\in\N$ there exist $k\in\N$
and elements $a^{(n)}_i\in J^{ni}$ for $i=1,\dots,k$,
such that $f^{k+1}+a^{(n)}_1f^k+\cdots+a^{(n)}_k=0$.
Therefore :
$$
|f^{k+1}|^*_J\leq
\max(|a^{(n)}_i|^*_J\cdot|f^{k+1-i}|^*_J~|~i=1,\dots,k)
$$
(lemma \ref{lem_normalized-Samuel}(ii)) whence
$(|f|^*_J)^i\leq|a^{(n)}_i|^*_J\leq\rho^{ni}$ for some
$i\leq k$, (lemma \ref{lem_normalized-Samuel}(iii)),
and consequently $|f|^*_J\leq\rho^n$ for every $n\in\N$,
whence (iii.a).
\end{proof}

\begin{corollary}\label{cor_ker-is-in-plus}
Let $\underline X:=(X,\cT_X,A^+_X)$ be any quasi-affinoid
scheme, and $(X^\wedge,\cT^\wedge_X,A^+_{X^\wedge})$ the completion
of\/ $\underline X$. Denote by $A^\wedge_X$ (resp. $A^{\wedge+}_X$)
the separated completion of $A_X:=\cO_{\!X}(X)$ (resp.
of $A^+_X$). Then we have :
\begin{enumerate}
\item
The kernel $\cJ$ of the natural map
$\pi:A^\wedge_X\to\cO_{\!X^\wedge}(X^\wedge)$ is contained in
$(A^\wedge_X)^{\circ\circ}$.
\item
$\pi^{-1}A^+_{X^\wedge}=A^{\wedge+}_X$.
\end{enumerate}
\end{corollary}
\begin{proof}(i): To begin with, we remark, quite generally :

\begin{claim}\label{cl_density-crit}
Let $A$ be any ring, $I\subset A$ a finitely generated
ideal, and set
$$
X_A:=\Spec\,A
\qquad
Z:=\Spec\,A/I
\qquad
U:=X_A\setminus Z.
$$
Then the following conditions are equivalent :
\begin{enumerate}
\alphaenu
\item
$U$ is dense in $X_A$.
\item
The support of the $A$-module $I^n/I^{n+1}$ is $Z$, for
every $n\in\N$. 
\end{enumerate}
\end{claim}
\begin{pfclaim} Since $U$ is constructible in $X_A$,
condition (a) is equivalent to :
\begin{itemize}
\item[(c)]
$I_\fq=A_\fq$ for every minimal prime ideal $\fq$ of $A$
\end{itemize}
(proposition \ref{prop_closed-under-spec}(iii)). Now,
if $A=0$ there is nothing to prove, so assume that $A\neq 0$.
Suppose that (c) holds. Clearly the support of $I^n/I^{n+1}$
lies in $Z$, so we consider any $\fp\in Z$, and we have
to show that $I^n_\fp/I^{n+1}_\fp\neq 0$ for every $n\in\N$.
However, if the latter fails for some $n\in\N$, Nakayama's
lemma implies that $I^n_\fp=0$; then $I^n_\fq=0$ for every
minimal prime ideal $\fq$ of $A$ contained in $\fp$,
contradicting (c).

Conversely, if (c) fails for some minimal prime ideal
$\fq$, it follows that $I_\fq$ is a nilpotent ideal of
$A_\fq$. Say that $I^n_\fq\neq 0$ and $I^{n+1}_\fq=0$ for
some $n\in\N$; then $\fq$ lies in the support of
$I^n/I^{n+1}$, hence $\fq\in Z$, and therefore $U$
cannot be dense in $X_A$.
\end{pfclaim}

Set $Y:=\Spec\,A_X$, $Y:=\Spec\,A^\wedge_X$, $X':=X\times_YY'$,
and let $I\subset A_X$ be any ideal such that
$A^{\circ\circ}_X\subset I$ and $Y\setminus X=\Spec\,A_X/I$.
Then $Y'\setminus X'=\Spec\,A^\wedge_X/I'$, where
$I':=IA^\wedge_X$, and since $I$ is open in $A_X$,
the natural maps $I^n/I^{n+1}\to I'^n/I'^{n+1}$ are
bijective for every $n\in\N$. Moreover, by assumption
$X$ is dense in $Y$; taking into account claim
\ref{cl_density-crit}, we deduce that $X'$ is dense
in $Y'$. Pick any finitely generated ideal $J$ of
$A^\wedge_X$ whose radical is $A_X^{\wedge\circ\circ}$;
in view of corollary \ref{cor_corcor}(ii,iii), we
see that
$$
\bigcap_{v\in\Spa\,(A^\wedge_X,A^{\wedge+}_X,X')}\Ker\,v=
\{a\in A^\wedge_X~|~|a|^*_J=0\}\subset(A^\wedge_X)^{\circ\circ}.
$$
But we have already noticed that $\cJ$ is contained
in the support of every $v\in\Spa(A^\wedge_X,A^{\wedge+}_X,X')$
(see the proof of lemma \ref{lem_invariance-by-loc-hens}(iii)),
whence (i).

(ii) follows immediately from (i) and corollary
\ref{cor_corcor}(i) and lemma \ref{lem_invariance-by-loc-hens}(iii)
(details left to the reader; more directly, one can also argue
as in the proof of theorem \ref{th_int-subrings-perfectoid}(ii)).
\end{proof}

\begin{proposition}\label{prop_crit-invertible}
Let $(X,\cT_X,A^+_X)$ be any topologically local quasi-affinoid
scheme, and $f\in\cO_{\!X}(X)$ any element. The following conditions
are equivalent :
\begin{enumerate}
\alphaenu
\item
$f\in\cO_{\!X}(X)^\times$.
\item
$v(f)\neq 0$ for every $v\in\Spa(X,\cT_X,A^+_X)$.
\end{enumerate}
\end{proposition}
\begin{proof}Clearly, it suffices to check that
(b)$\Rightarrow$(a). Hence, suppose that (b) holds,
and set
$$
A_X:=\cO_{\!X}(X)
\qquad
Y:=\Spec\,A_X
\qquad
Y^{\circ\circ}:=\Spec\,A_X/A^{\circ\circ}_X\cdot A_X.
$$
Let $x\in X$ be any point; we need to show that the
image $f(x)$ of $f$ does not vanish in $\kappa(x)$.
If $x\in Y^{\circ\circ}$, let $v$ be the trivial valuation
with support equal to $x$; clearly $v\in\Spa(X,\cT_X,A^+_X)$,
so (b) says that $v(f)\neq 0$, whence the contention,
in this case.

Next, if $x\in X\setminus Y^{\circ\circ}$, we proceed as in
the proof of proposition \ref{prop_f-adic-Spectra}(iii.b) :
we pick first a minimal specialization $y$ of $x$ in
$X\setminus Y^{\circ\circ}$, then a maximal point $\fp$ of the
non-empty closed subset $\{y\}^c\cap Y^{\circ\circ}$ of $Y$
(where $\{y\}^c$ denotes the topological closure of $\{y\}$
in $Y$). The image $C$ of $A_{X,\fp}$ in $\kappa(y)$ is a
one-dimensional local domain, so we may find a valuation
ring $V$ of $\kappa(y)$ with value group of rank one, that
dominates $C$ (claim \ref{cl_one-dim-dominates}). Let
$v\in\Spv\,A_X$ be the valuation corresponding to $V$;
arguing as in {\em loc. cit.} we see that
$v\in\Spa(X,\cT_X,A^+_X)$, hence $v(f)\neq 0$, so
$f(y)\neq 0$ and finally $f(x)\neq 0$ as well.
\end{proof}

\subsection{Adic spaces}\label{subsec_adic-spaces}
In this section we shall endow the adic spectra introduced
in section \ref{sec_affinoid-rings} with certain natural
presheaves of topological rings, and we shall show that,
under suitable conditions, these presheaves are sheaves of
rings or of topological rings. The first step is the following
universal construction :

\sset\subsubsection{}\label{subsec_universal-property}
Let $\underline A:=(A,A^+,U)$ be any quasi-affinoid ring,
$\Lambda$ a (small) set, $S:=(s_\lambda~|~\lambda\in\Lambda)$
a system of elements of $A$, and $T:=(T_\lambda~|~\lambda\in\Lambda)$
a family of subsets of $A$ such that $T_\lambda$ generates an
open ideal of $A$, for every $\lambda\in\Lambda$. We consider
the f-adic ring $A[X_\lambda~|~\lambda\in\Lambda]_T$
provided by proposition \ref{prop_polynomial-top-rings}(ii),
and its ideal $J$ generated by the system
$(1-s_\lambda X_\lambda~|~\lambda\in\Lambda)$. We set
$$
A\Bigl(\frac{T}{S}\Bigr):=A[X_\bullet]_T/J
$$
and we endow this ring with the quotient topology
induced by $A[X_\bullet]_T$ via the natural projection.
In other words, $A(\frac{T}{S})$ is an f-adic topological
ring whose underlying $A$-algebra is naturally identified
with the localization $A[s_\lambda^{-1}~|~\lambda\in\Lambda]$.
From proposition \ref{prop_polynomial-top-rings}(iii,iv) we
see that -- under this identification -- the subset
$\{t/s_\lambda~|~\lambda\in\Lambda,\ t\in T_\lambda\}$
is power bounded in $A(\frac{T}{S})$. Moreover,
the localization map $h:A\to A(\frac{T}{S})$ is f-adic,
and enjoys the following universal property.
For every f-adic ring $B$ and every continuous map
$f:A\to B$ such that :
\begin{itemize}
\item
$f(s_\lambda)\in B^\times$ for every $\lambda\in\Lambda$
\item
the subset $\{f(t)/f(s_\lambda)~|~
\lambda\in\Lambda,\ t\in T_\lambda\}$ is power bounded
in $B$
\end{itemize}
there exists a unique continuous map
$g:A(\frac{T}{S})\to B$ such that $g\circ h=f$. Furthermore,
if $A_0\subset A$ is a subring of definition,
$A_0[\frac{t}{s_\lambda}~|~\lambda\in\Lambda,\ t\in T_\lambda]$
is a subring of definition of $A(\frac{T}{S})$.

Next, notice that
$C:=A^+[\frac{t}{s_\lambda}~|~\lambda\in\Lambda,\ t\in T_\lambda]$
is an open subring of $A(\frac{T}{S})^\circ$, and denote
by $C'$ the integral closure of $C$ in $A(\frac{T}{S})$.
Moreover, set $X:=\Spec\,A[s_\lambda^{-1}|\lambda\in\Lambda]$,
and notice that $U':=U\cap X$ contains the analytic locus
of $X$, since $h$ is f-adic (lemma \ref{lem_deja-vu}(iv)).
We obtain therefore a quasi-affinoid ring
$$
\underline A\Bigl(\frac{T}{S}\Bigr):=
\Bigl(A\Bigl(\frac{T}{S}\Bigr),C',U'\Bigr)
$$
and the localization map yields a natural f-adic morphism
of quasi-affinoid rings
\set\begin{equation}\label{eq_univ-map-aff}
\underline A\to\underline A\Bigl(\frac{T}{S}\Bigr).
\end{equation}
In case $\Lambda=\{\lambda\}$ has only one element,
$T_\lambda=\{f_0,\dots,f_n\}$ is a finite subset of $A$ which
generates an open ideal, and $S=\{s\}$ for some $s\in A$,
we also denote this quasi-affinoid ring by :
$$
\underline A\Bigl(\frac{f_0,\dots,f_n}{s}\Bigr).
$$

\sset\subsubsection{}\label{subsec_represent-rationals}
We consider now the Yoneda embedding (see \eqref{subsec_yoneda})
$$
h:\mathsf{q.Afd.Sch}_\loc\to
\mathsf{q.Afd.Sch}^\wedge_\loc
\qquad
\underline X\mapsto h_{\underline X}
$$
(notation of remark \ref{rem_loc-hens-cplete}(v)).
Let $\underline X$ be any topologically local quasi-affinoid
scheme, and $U\subset\Spa\,\underline X$ any open subset.
We attach to $U$ the sub-presheaf $h_U\subset h_{\underline X}$
given by the rule :
$$
h_U(\underline Y):=
\{\phi\in h_{\underline X}(\underline Y)~|~
  \Img\,\Spa(\phi)\subset U\}
\qquad
\text{for every $\underline Y\in\Ob(\mathsf{q.Afd.Sch}_\loc)$}.
$$

\begin{definition}
(i)\ \
In the situation of \eqref{subsec_represent-rationals},
we say that $U$ is a {\em quasi-affinoid open subset\/}
of $\Spa\,\underline X$, if the presheaf $h_U$ is
representable by some $\underline Y\in\Ob(\mathsf{q.Afd.Sch}_\loc)$.

(ii)\ \
We say that $U$ is an {\em affinoid open subset\/} of
$\Spa\,\underline X$, if $\underline Y$ as in (i) is an
affinoid scheme.
\end{definition}

\begin{theorem}\label{th_represent-rational}
Every rational subset of\/ $\Spa\,\underline X$ is a
quasi-affinoid open subset.
\end{theorem}
\begin{proof} Let $R$ be such a rational subset, and say that
$$
\sGamma(\underline X)=(A_X,A_X^+,X)
\qquad\text{and}\qquad
R=R_{A_X}\Bigl(\frac{f_1}{f_0},\dots,\frac{f_n}{f_0}\Bigr)
\cap\Spa\,\underline X
$$
for a sequence $(f_0,f_1,\dots,f_n)$ of elements of $A_X$
that generates an open ideal. We consider the quasi-affinoid
rings
$$
\underline A_R:=
\sGamma(\underline X)\Bigl(\frac{f_0,\dots,f_n}{f_0}\Bigr)
$$
defined as in \eqref{subsec_universal-property}. Also, we let
$$
i_R:\sGamma(\underline X)\to\underline A_{R,\loc}
$$
be the composition of the natural morphism
$\sGamma(\underline X)\to\underline A_R$ of
\eqref{eq_univ-map-aff} with the topological localization
$\underline A_R\to\underline A_{R,\loc}$.
Now, let $\underline Y$ be any topologically local
quasi-affinoid scheme, and $\phi\in h_R(\underline Y)$
any element. Set $(A_Y,A^+_Y,Y):=\sGamma(\underline Y)$,
and let $\phi^\flat:A_X\to A_Y$ be the continuous map
associated with $\phi$.

\begin{claim}\label{cl_good-situation}
$\phi^\flat(f_0)\in A_Y^\times$ and
$\phi^\flat(f_i)/\phi^\flat(f_0)\in A_Y^+$ for every $i=1,\dots,n$.
\end{claim}
\begin{pfclaim} Indeed, for every $v\in\Spa\,\underline Y$
we have by assumption $v\circ\phi^\flat\in R$, so that
$v(\phi^\flat(f_0))\neq 0$ and
$v(\phi^\flat(f_i))\leq v(\phi^\flat(f_0))$ for $i=1,\dots,n$.
Then proposition \ref{prop_crit-invertible}
implies that $\phi^\flat(f_0)\in A_Y^\times$, so we may write
$v(\phi^\flat(f_i)/\phi^\flat(f_0))\leq 1$ for every $i=1,\dots,n$
and every $v\in\Spa\,\underline Y$. To conclude, it
then suffices to invoke corollary \ref{cor_corcor}(i).
\end{pfclaim}

From claim \ref{cl_good-situation}, remark
\ref{rem_something-on-bdd}(ii) and the discussion of
\eqref{subsec_universal-property} we see that $\phi^\flat$
factors uniquely through a continuous ring homomorphism
$g:A_R\to A_Y$, and $g(f_i/f_0)\in A^+_Y$ for every
$i=1,\dots,n$. Thus, $g$ defines a morphism of quasi-affinoid
rings
$$
g:\underline A_R\to\sGamma(\underline Y)
$$
and summing up, we easily conclude that $\phi$ factors
uniquely through $\sSpec(i_R)$ and the morphism of
topologically local quasi-affinoid schemes
$$
\psi:=\sSpec(g_\loc):\underline Y\to
\underline R:=\sSpec\,(\underline A_{R,\loc}).
$$
Lastly, taking into account lemma
\ref{lem_invariance-by-loc-hens}(i) it is easily seen that
$\Spa(i_R)$ is an injective map inducing a homeomorphism
$\Spa\,\underline R\isom R$, and the theorem follows.
\end{proof}

\sset\subsubsection{}\label{subsec_quasi-affoid-subsets}
In the situation of \eqref{subsec_represent-rationals},
let $U\subset\Spa\,\underline X$ be a quasi-affinoid
open subset, and $\underline Y$ any topologically local
quasi-affinoid scheme that represents $h_U$. Then the
morphism $\one_{\underline Y}$ corresponds to a well defined
section $\phi_{\underline Y/\underline X}\in h_U(\underline Y)$,
{\em i.e.} to a morphism
$\phi_{\underline Y/\underline X}:\underline Y\to\underline X$ of
quasi-affinoid schemes with
$\Img(\Spa\,\phi_{\underline Y/\underline X})\subset U$.

\begin{corollary}\label{cor_rational-represent}
With the notation of \eqref{subsec_quasi-affoid-subsets},
the following holds :
\begin{enumerate}
\item
The map
$\Spa\,\phi_{\underline Y/\underline X}:
\Spa\,\underline Y\to\Spa\,\underline X$
induces a homeomorphism $\Spa\,\underline Y\isom U$.
\item
Especially, every quasi-affinoid open subset is quasi-compact.
\item
More precisely, $\Spa\,\phi_{\underline Y/\underline X}$ induces
homeomorphisms
$$
(\Spa\,\underline Y)_\mathrm{a}\isom U\cap\Spa(\underline X)_\mathrm{a}
\qquad
(\Spa\,\underline Y)_\mathrm{na}\isom U\cap\Spa(\underline X)_\mathrm{na}.
$$
\item
The morphism $\phi_{\underline Y/\underline X}$ is f-adic.
\end{enumerate}
\end{corollary}
\begin{proof}(i): Let $(R_i~|~i\in I)$ be a system of rational
subsets of $\Spa\,\underline X$ such that $U=\bigcup_{i\in I}R_i$,
and for every $i,j\in I$ set $R_{ij}:=R_i\cap R_j$ and choose
topologically local quasi-affinoid schemes $\underline Z{}_i$
and $\underline Z{}_{ij}$ that represent $h_{R_i}$ and $h_{R_{ij}}$.
The inclusions $R_i\subset U$ and $R_{ij}\subset U$ are then
represented by morphisms $\psi_i:\underline Z{}_i\to\underline Y$
and respectively $\psi_{ij}:\underline Z{}_{ij}\to\underline Y$
of quasi-affinoid schemes such that
$\phi_{\underline Z{}_i/\underline X}=\phi_{\underline Y/\underline X}\circ\psi_i$
and $\phi_{\underline Z{}_{ij}/\underline X}=
\phi_{\underline Y/\underline X}\circ\psi_{ij}$, for every $i,j\in I$.
For every $i\in I$ set
$U_i:=\Img\,\Spa\,\psi_i\subset\Spa\,\underline Y$; it was
observed in the proof of theorem \ref{th_represent-rational}
that $\Spa\,\phi_{\underline Z{}_i/\underline X}$ induces a homeomorphism
$\Spa\,\underline Z{}_i\isom R_i$; it follows easily that
$\Spa\,\psi_i$ and $\Spa\,\phi_{\underline Y/\underline X}$ induce
homeomorphisms $\Spa\,\underline Z{}_i\isom U_i$ and
respectively $U_i\isom R_i$. It is also easily seen that
$U_i\cap U_j=\Img\,\Spa\,\psi_{ij}$, so the map
$\Spa\,\phi_{\underline Y/\underline X}$ restricts to a homeomorphism
$U':=\bigcup_{i\in I}U_i\isom U$. To conclude, it suffices
to check that $U'=\Spa\,\underline Y$. However, suppose
that $v\in\Spa\,\underline Y\setminus U'$; we may find
a rational subset $R'$ of $\Spa\,\underline Y$ that
contains $v$, and we may assume that
$\Spa\,\phi_{\underline Y/\underline X}(R')\subset R_i$ for some
$i\in I$. Pick a topologically local quasi-affinoid scheme
$\underline Y'$ that represents the subsheaf $h_{R'}$
of $h_{\underline Y}$, and let $\phi':\underline Y'\to\underline Y$
be the morphism of quasi-affinoid schemes that represents
the inclusion $R'\subset\Spa\,\underline Y$. Then
$\phi_{\underline Y/\underline X}\circ\phi':\underline Y'\to\underline X$
lies in $h_{R_i}(\underline Y')$, so it corresponds to a
unique morphism $\beta:\underline Y'\to\underline Z{}_i$
such that
$\phi_{\underline Y/\underline X}\circ\phi'=
\phi_{\underline Z{}_i/\underline X}\circ\beta=
\phi_{\underline Y/\underline X}\circ\psi_i\circ\beta$. Since
$\phi_{\underline Y/\underline X}$ induces an injective morphism
$h_{\underline Y}\to h_{\underline X}$ (whose image is $h_U$),
we deduce that $\phi'=\psi_i\circ\beta$. Especially, $U_i$
contains the image of $\Spa\,\phi'$; but we know that the
latter coincides with $R'$, a contradiction.

(ii) is an immediate consequence of (i).

(iii): By virtue of remark \ref{rem_Cont-A}(vii), it
suffices to show that $\Spa\,\phi_{\underline Y/\underline X}$ maps
$(\Spa\,\underline Y)_\mathrm{a}$ into $(\Spa\,X)_\mathrm{a}$.
Thus, let $v\in(\Spa\,\underline Y)$, pick $i\in I$
such that $w:=\Spa\,\phi_{\underline Y/\underline X}(v)\in R_i$,
and let $u:=(\Spa\,\phi_{\underline Z{}_i/\underline X})^{-1}(w)$.
Suppose that $w\in(\Spa\,X)_\mathrm{na}$; by inspecting
the proof of theorem \ref{th_represent-rational}, it is
easily seen that $\phi_{\underline Z{}_i/\underline X}$ is f-adic,
so $u\in\Spa(\underline Z{}_i)_\mathrm{na}$, by corollary
\ref{cor_Cont-spectral}(ii). By the foregoing, we also
know that $\Spa\,\psi_i(u)=v$, and then
$v\in(\Spa\,\underline Y)_\mathrm{na}$, again by remark
\ref{rem_Cont-A}(vii), whence the assertion.

(iv) follows immediately from (iii) and proposition
\ref{prop_f-adic-Spectra}(iii.b).
\end{proof}

\begin{remark}\label{rem_yoneda-rationals}
(i)\ \
In the situation of \eqref{subsec_represent-rationals}, we
can also consider the Yoneda imbeddings for the categories
of topologically henselian and of complete, separated
quasi-affinoid schemes
$$
h':\mathsf{q.Afd.Sch}_\mathrm{hens}\to
\mathsf{q.Afd.Sch}^\wedge_\mathrm{hens}
\qquad
h'':\mathsf{q.Afd.Sch}_\mathrm{comp}\to
\mathsf{q.Afd.Sch}^\wedge_\mathrm{comp}
$$
and if $\underline X$ is topologically henselian
(resp. complete and separated) we set
$$
h'_U(\underline Y'):=h_U(\underline Y')\subset
h'_{\underline X}(\underline Y')
\qquad
\text{(resp.\  $h''_U(\underline Y''):=h_U(\underline Y'')
\subset h''_{\underline X}(\underline Y'')$\ )}
$$
for every topologically henselian quasi-affinoid scheme
$\underline Y'$ and every complete and separated quasi-affinoid
scheme $\underline Y''$. Suppose that $U$ is quasi-affinoid,
and let $\phi_{\underline Y/\underline X}:\underline Y\to\underline X$
be a morphism of topologically local quasi-affinoid schemes
representing the sub-presheaf $h_U\subset h_{\underline X}$,
as in \eqref{subsec_quasi-affoid-subsets}; from remark
\ref{rem_loc-hens-cplete}(v) and lemma
\ref{lem_invariance-by-loc-hens}(i) we see that if $\underline X$
is topologically henselian, the topological henselization
$\phi^\he_{\underline Y/\underline X}:\underline Y^\he\to\underline X$
of $\phi_{\underline Y/\underline X}$ represents $h'_U$. Likewise,
taking into account proposition \ref{prop_complete-qaff-sch}
and lemma \ref{lem_invariance-by-loc-hens}(iii), we conclude
that if $\underline X$ is complete and separated, the completion
$\phi^\wedge_{\underline Y/\underline X}:\underline Y^\wedge\to\underline X$
of $\phi_{\underline Y/\underline X}$ represents $h''_U$.

(ii)\ \
Moreover, we may consider the Yoneda imbedding for the
opposite of the category of topologically local affinoid
rings :
$$
h''':\mathsf{Afd.Ring}^o_\loc\to
\mathsf{Afd.Ring}^{o\wedge}_\loc
$$
and if $\underline A:=(A,A^+)$ is any topologically
local affinoid ring, and $U\subset\Spa\,\underline A$
any open subset, we may form again the sub-presheaf
$h'''_U:\subset h'''_{\underline A}$ by the rule :
$h'''_U(\underline B):=h_U(\sSpec\,\underline B)$ for
every topologically local affinoid ring $\underline B$.
We may regard $U$ as an open subset of
$\Spa(\sSpec\,\underline A)$, and we say that $U$ is
a quasi-affinoid open subset of $\Spa\,\underline A$ if it
is such, as an open subset of $\Spa\,\sSpec\,\underline A$.
Since $\sSpec$ is fully faithful on the full subcategory
$\mathsf{Afd.Ring}$, it follows that if $U$ is quasi-affinoid,
$h'''_U$ is representable : namely, a representing
affinoid ring is $\sGamma(\underline Y)$, where $\underline Y$
is any topologically local quasi-affinoid scheme that represents
the presheaf $h_U$. Lastly, we may in the same way consider
variants for topologically henselian (resp. complete
and separated) affinoid rings : the reader may spell out
the details. (However, the method does not produce representing
objects on the whole category of quasi-affinoid rings,
essentially because there is -- in this generality -- no
counterpart of proposition \ref{prop_crit-invertible}
for quasi-affinoid rings).

(iii)\ \
It follows from lemma \ref{lem_f-adic-pushout} and corollary
\ref{cor_rational-represent}(iv) that the intersection of two
quasi-affinoid open subsets $U,U'\subset\Spa\,\underline X$
is quasi-affinoid. Indeed, if $\underline Y$ and
$\underline Y'$ are topologically local quasi-affinoid
schemes that represent $h_U$ and respectively $h_{U'}$,
it is easily seen that
$$
(\underline Y\times_{\underline X}\underline Y')_\loc
$$
represents $h_U\cap h_{U'}=h_{U\cap U'}$ (see example
\ref{ex_fibre-prod-in-qaff-sch}).

(iv)\ \
In the same vein, if $f:\underline X'\to\underline X$ is any
f-adic morphism of topologically local quasi-affinoid schemes,
and $U\subset\Spa\,\underline X$ is any quasi-affinoid open
subset, then $U':=(\Spa\,f)^{-1}U$ is a quasi-affinoid open
subset of $\Spa\,\underline X'$. Indeed, if $\underline Y$
represents $h_U$, then $h_{U'}$ is represented by 
$(\underline Y\times_{\underline X}\underline X')_\loc$, where
the fibre product is given as in example
\ref{ex_fibre-prod-in-qaff-sch}.

(v)\ \
Lastly, let $U$ be a quasi-affinoid open subset of
$\Spa\,\underline X$ such that the inclusion
$h_U\subset h_{\underline X}$ is represented by the morphism
of topologically local quasi-affinoid schemes
$\phi_{\underline Y/\underline X}:\underline Y\to\underline X$,
and $V$ a quasi-affinoid open subset of
$\Spa\,\underline Y$ such that the inclusion
$h_V\subset h_{\underline Y}$ is represented by the morphism
of topologically local affinoid schemes
$\phi_{\underline Z/\underline Y}:\underline Z\to\underline Y$,
and set $U':=\Spa\,\phi_{\underline Y/\underline X}(V)\subset U$.
Then $\phi_{\underline Y/\underline X}$ induces an isomorphism
$h_{\underline Y}\isom h_U$ that identifies $h_V$ with the
sub-presheaf $h_{U'}$ of $h_U$. Especially, $U'$ is a
quasi-affinoid open subset of $\Spa\,X$, and the inclusion
$h_{U'}\subset h_{\underline X}$ is represented by
$\phi_{\underline Y/\underline X}\circ\phi_{\underline Z/\underline Y}$.
\end{remark}

\begin{lemma}\label{lem_rat-in-rat}
In the situation of remark {\em\ref{rem_yoneda-rationals}(v)},
suppose furthermore that $U$ is a rational subset of\/
$\Spa\,\underline X$ and $V$ is a rational subset of\/
$\Spa\,\underline Y$. Then $V$ is a rational subset of\/
$\Spa\,\underline X$.
\end{lemma}
\begin{proof} Say that $\sGamma(\underline X)=(A,A^+,X)$,
and $U=R_A(\frac{f_1}{f_0},\dots,\frac{f_n}{f_0})\cap
\Spa\,\underline X$, with $(f_0,\dots,f_n)$ a sequence of
elements of $A$ that generates an open ideal. Set
$(A_U,A^+_U):=
\underline A(\frac{f_1}{f_0},\dots,\frac{f_n}{f_0})$, where
$\underline A:=(A,A^+)$, so that $A_U=A[f_0^{-1}]$. By
inspecting the proof of theorem \ref{th_represent-rational},
we get an isomorphism of quasi-affinoid schemes
$$
\underline Y\isom\sSpec\,(A_U,A^+_U,X\cap\Spec\,A_U)_\loc
$$
Since $X$ is quasi-compact, and $\Spec\,A_U$ is an open
affine subset of $\Spec\,A$, we have
$\cO_{\!X}(X\cap\Spec\,A_U)=\cO_{\!X}(X)[f_0^{-1}]=A_U$
(\cite[Ch.I, Prop.9.2.1]{EGAI}). Taking into account
\eqref{eq_wrong-order}, we deduce that
$\sGamma(\underline Y)=(B,B^+,X\cap\Spec\,A_U)$, where
$B:=A_{U,\loc}$ and $B^+:=A^+_{U,\loc}$. Consequently, we have
$V=R_B(\frac{g_1}{g_0},\dots,\frac{g_m}{g_0})\cap
\Spa\,\underline Y$ for a sequence $g_\bullet:=(g_0,\dots,g_m)$
of elements of $B$ that generates an open ideal $J$. Since
$B$ is a localization of $A$, we may also find $b\in B^\times$,
such that $bg_i\in A':=\Img(A\to B)$; after replacing
$g_\bullet$ by the sequence $(bg_0,\dots,bg_n)$, we may
then assume that $g_0,\dots,g_m\in A'$, and for every
$i=0,\dots,m$ we pick $h_i\in A$ whose image in $B$
agrees with $g_i$. Moreover, let $A_0\subset A$ be a ring
of definition, and $I\subset A_0$ a finitely generated
ideal of adic definition; according to claim
\ref{cl_efficient}(i), there exists $k\in\N$ such that
$v(a)<v(h_0)$ for every $a\in I^k$ and every $v\in V$.
Pick a finite system $x_1,\dots,x_r$ of generators for
$I^k$ and set $h_{i+m}:=x_i$ for every $i=1,\dots,r$; by
construction, we easily see that
$$
V=U\cap R
\qquad
\text{where
$R:=R_A\Bigl(\frac{h_1}{h_0},\dots,\frac{h_{m+r}}{h_0}\Bigr)$}
$$
and the system $(h_i~|~i=1,\dots,m+r)$ generates an open
ideal of $A$, so $R$ is a rational subset of
$\Spa\,\underline X$, and finally, the same holds for $V$.
\end{proof}

\sset\subsubsection{}\label{subsec_rational-site}
Let $\underline X:=(X,\cT_X,A^+_X)$ be any topologically
local quasi-affinoid scheme; in view of remark
\ref{rem_yoneda-rationals}(iii) we may associate with
$\underline X$ a site
$$
(\cQ(\underline X),J_\cQ)
$$
where $\cQ(\underline X)$ is the category of all quasi-affinoid
open subsets of $\Spa\,\underline X$; for every
$U,U'\in\Ob(\cQ(\underline X))$, the set
$\Hom_{\cQ(\underline X)}(U,U')$ contains exactly one morphism if
$U\subset U'$, and is empty otherwise. For every such $U$, the
sieves covering $U$ for the topology $J_\cQ$ are precisely those
generated by the families $(U_\lambda~|~\lambda\in\Lambda)$ of
quasi-affinoid open subsets of $\Spa\,\underline X$ such that
$\bigcup_{\lambda\in\Lambda}U_\lambda=U$.
We consider nine presheaves of topological $\cO_{\!X}(X)$-algebras
on $\cQ(\underline X)$, related by natural morphisms of presheaves :
$$
\xymatrix{ \cO^{\loc+}_{\Spa\,\underline X} \ar@{^{(}->}[r] \ar[d] &
\cO^{\loc\,\circ}_{\Spa\,\underline X} \ar@{^{(}->}[r] \ar[d] &
\cO^\loc_{\Spa\,\underline X} \ar[d] \\
\cO^{\he+}_{\Spa\,\underline X} \ar@{^{(}->}[r] \ar[d] &
\cO^{\he\,\circ}_{\Spa\,\underline X} \ar@{^{(}->}[r] \ar[d] &
\cO^\he_{\Spa\,\underline X} \ar[d] \\
\cO^{\wedge+}_{\Spa\,\underline X} \ar@{^{(}->}[r] &
\cO^{\wedge\circ}_{\Spa\,\underline X} \ar@{^{(}->}[r] &
\cO^\wedge_{\Spa\,\underline X}
}$$
and constructed as follows. By theorem \ref{th_represent-rational},
for every $U\in\Ob(\cQ(\underline X))$ the sub-presheaf
$h_U$ of $h_{\underline X}$  is representable (notation
of \eqref{subsec_represent-rationals}); we choose a
topologically local quasi-affinoid scheme
$\underline X{}_U:=(X_U,\cT_U,A^+_U)$ that represents this
presheaf, and we set
$$
\cO^\loc_{\Spa\,\underline X}(U):=\cO_{\!X_U}(X_U)
\qquad
\cO^{\loc\,\circ}_{\Spa\,\underline X}(U):=\cO_{\!X_U}(X_U)^\circ
\qquad
\cO^{\loc+}_{\Spa\,\underline X}(U):=A^+_U.
$$
Any inclusion $U\subset U'$ of quasi-affinoid open subsets
induces a morphism $h_U\to h_{U'}$ of presheaves, which is
represented by a unique morphism
$\underline X{}_U\to\underline X{}_{U'}$ of topologically
local quasi-affinoid schemes (see remark
\ref{rem_represent-morph}(i)), whence well defined continuous
homomorphisms
$$
\cO^{\loc+}_{\Spa\,\underline X}(U')\to\cO^{\loc+}_{\Spa\,\underline X}(U)
\qquad
\cO^\loc_{\Spa\,\underline X}(U')\to\cO^\loc_{\Spa\,\underline X}(U)
$$
of topological $\cO_{\!X}(X)$-algebras, and clearly
the rules $U\mapsto\cO^{\loc+}_{\Spa\,\underline X}(U)$ and
$U\mapsto\cO^\loc_{\Spa\,\underline X}(U)$ yield two functors
from $\cQ(\underline X)^o$ to the category of topologically
local f-adic $\cO_{\!X}(X)$-algebras. Likewise, in view of
corollary \ref{cor_rational-represent}(iv) and lemma
\ref{lem_f-adics}(iii.a) we also get a continuous map
$\cO^{\loc\,\circ}_{\Spa\,\underline X}(U')\to
\cO^{\loc\,\circ}_{\Spa\,\underline X}(U)$, whence the presheaf
$\cO^{\loc\,\circ}_{\Spa\,\underline X}$.

Similarly, let $\underline X{}^\he_U:=(X^\he_U,\cT^\he_U,A^{\he+}_U)$
(resp. $\underline X{}^\wedge_U:=(X^\wedge_U,\cT^\wedge_U,A^{\wedge+}_U)$)
be a topologically henselian (resp. complete and separated)
quasi-affinoid scheme representing the sub-presheaf $h'_U$
of $h'_{\underline X^\he}$ (resp. $h''_U$ of $h''_{\underline X^\wedge}$)
as in remark \ref{rem_yoneda-rationals}(i); we define
$$
\begin{aligned}
\cO^\he_{\Spa\,\underline X}(U):=\,&
\cO_{\!X^\he_U}(X^\he_U)
\qquad &
\cO^{\he\,\circ}_{\Spa\,\underline X}(U):=\,&
\cO_{\!X^\he_U}(X^\he_U)^\circ
\qquad &
\cO^{\he+}_{\Spa\,\underline X}(U):=\,&
A^{\he+}_U \\
\cO^\wedge_{\Spa\,\underline X}(U):=\,&
\cO_{\!X^\wedge_U}(X^\wedge_U)
\qquad &
\cO^{\wedge\circ}_{\Spa\,\underline X}(U):=\,&
\cO_{\!X^\wedge_U}(X^\wedge_U)^\circ
\qquad &
\cO^{\wedge+}_{\Spa\,\underline X}(U):=\,&
A^{\wedge+}_U
\end{aligned}
$$
and arguing as in the foregoing, it is easily seen that these
rules yield well defined presheaves of topologically henselian
(resp. complete and separated) f-adic $\cO_{\!X}(X)$-algebras.
All these presheaves depend on the choices of representing
quasi-affinoid schemes, but for any two such sets of choices
there exists a unique isomorphism of presheaves of topological
$\cO_{\!X}(X)$-algebras between the corresponding presheaves.
We notice as well that the rule
$$
U\mapsto\cO^\loc_{\Spa\underline X}(U)^{\circ\circ}
$$
defines a sub-presheaf $\cO^{\loc\,\circ\circ}_{\Spa\,\underline X}
\subset\cO^\loc_{\Spa\,\underline X}$ that is a presheaf of ideals
in both $\cO^{\loc+}_{\Spa\,\underline X}$
and $\cO^{\loc\,\circ}_{\Spa\,\underline X}$. Likewise we get
presheaves of topologically nilpotent sections
$\cO^{\he\,\circ\circ}_{\Spa\,\underline X}\subset
\cO^\he_{\Spa\,\underline X}$ (resp.
$\cO^{\wedge\,\circ\circ}_{\Spa\,\underline X}\subset
\cO^\wedge_{\Spa\,\underline X}$) that are presheaves of ideals
in $\cO^{\he\,\circ}_{\Spa\,\underline X}$ and
$\cO^{\he+}_{\Spa\,\underline X}$ (resp. in
$\cO^{\wedge\,\circ}_{\Spa\,\underline X}$ and
$\cO^{\wedge+}_{\Spa\,\underline X}$).

\sset\subsubsection{}\label{subsec_valuation-on-stalks}
With the notation of \eqref{subsec_rational-site}, let
$x\in\Spa\,\underline X$ be any point; recall that $x$
is the equivalence class of a continuous valuation
$v_x:A:=\cO_{\!X}(X)\to\Gamma_x$; if $R\subset\Spa\,\underline X$
is any rational subset containing $x$, the valuation $v_x$
extends uniquely to a continuous valuation $v^\wedge_{x,R}$
on $\cO^\wedge_{\Spa\,\underline X}(R)$ with value group $\Gamma_x$
(claim \ref{cl_complete-cont-val}(i)). Then, by restricting
along the natural maps
$\cO^\loc_{\Spa\,\underline X}(R)\to\cO^\he_{\Spa\,\underline X}(R)\to
\cO^\wedge_{\Spa\,\underline X}(R)$ we get corresponding valuations
$v^\he_{x,R}$, $v^\loc_{x,R}$, and -- after taking colimits --
unique valuations on the stalks of these presheaves
\set\begin{equation}\label{eq_vals-on-stalks}
{\diagram \cO^\loc_{\Spa\,\underline X,x} \ar[r]
\ar@/_1pc/[rd]_{|\cdot|^\loc_x} & \cO^\he_{\Spa\,\underline X,x} \ar[r]
\ar[d]^{|\cdot|^\he_x} &
\cO^\wedge_{\Spa\,\underline X,x} \ar@/^1pc/[ld]^{|\cdot|^\wedge_x} \\
& \Gamma_{\!x}.
\enddiagram}
\end{equation}
We denote by $\kappa(x^\loc)$, $\kappa(x^\he)$ and $\kappa(x^\wedge)$
the residue fields of $|\cdot|^\loc_x$, $|\cdot|^\he_x$ and
respectively $|\cdot|^\wedge_x$, and recall that the valuations
of \eqref{eq_vals-on-stalks} induce natural residual valuations
on these fields : see remark \ref{rem_semi-norm}(v). There follow
natural inclusion of valued fields
$$
\kappa(x^\loc)\to\kappa(x^\he)\to\kappa(x^\wedge)
$$
and taking into account claim \ref{cl_complete-cont-val}(ii)
we see that the images of $\kappa(x^\loc)$ and of $\kappa(x^\he)$
in $\kappa(x^\wedge)$ are both dense for the valuation topology
of the latter field. Especially, these inclusions induce natural
identifications of the completions of these fields for their
respective valuation topologies (theorem
\ref{th_complete-top-grps}(iii)), and we denote this common
completion by
$$
(\kappa(x)^\wedge,|\cdot|_x^\wedge)
$$
which is a valued field with value group $\Gamma_{\!x}$, by
proposition \ref{prop_stays-valuation}(iii,v).

\begin{lemma}\label{lem_stalks-are-local}
In the situation of \eqref{subsec_valuation-on-stalks}, we have :
\begin{enumerate}
\item
The stalks $\cO^\wedge_{\Spa\,\underline X,x}$, $\cO^\he_{\Spa\,\underline X,x}$
and $\cO^\loc_{\Spa\,\underline X,x}$ are local rings, and
$\kappa(x^\loc)$, $\kappa(x^\he)$, $\kappa(x^\wedge)$ are
their respective residue fields.
\item
The stalks $\cO^{\wedge+}_{\Spa\,\underline X,x}$,
$\cO^{\he+}_{\Spa\,\underline X,x}$ and $\cO^{\loc+}_{\Spa\,\underline X,x}$
are local rings, and the induced diagrams
$$
\xymatrix{ \cO^{\loc+}_{\Spa\,\underline X,x} \ar[r] \ar[d] &
\cO^\loc_{\Spa\,\underline X,x} \ar[d]^{\pi^\loc_x} &
\cO^{\he+}_{\Spa\,\underline X,x} \ar[r] \ar[d] &
\cO^\he_{\Spa\,\underline X,x} \ar[d]^{\pi^\he_x} &
\cO^{\wedge+}_{\Spa\,\underline X,x} \ar[r] \ar[d] &
\cO^\wedge_{\Spa\,\underline X,x} \ar[d]^{\pi^\wedge_x} \\
\kappa(x^\loc)^+ \ar[r] & \kappa(x^\loc) &
\kappa(x^\he)^+ \ar[r] & \kappa(x^\he) &
\kappa(x^\wedge)^+ \ar[r] & \kappa(x^\wedge)
}$$
are cartesian (here $\kappa(x^\loc)^+$ is the valuation ring
of the residual valuation on $\kappa(x^\loc)$, and likewise
for $\kappa(x^\he)^+$ and $\kappa(x^\wedge)^+$).
\item
The stalk $\cO^{\loc\,\circ\circ}_{\Spa\,\underline X,x}$ (resp.
$\cO^{\he\,\circ\circ}_{\Spa\,\underline X,x}$, resp.
$\cO^{\wedge\,\circ\circ}_{\Spa\,\underline X,x}$) is the radical
of the ideal $A^{\circ\circ}\cdot\cO^{\loc+}_{\Spa\,\underline X,x}$ of
$\cO^{\loc+}_{\Spa\,\underline X,x}$ (resp.
$A^{\circ\circ}\cdot\cO^{\he+}_{\Spa\,\underline X,x}$ of
$\cO^{\he+}_{\Spa\,\underline X,x}$, resp.
$A^{\circ\circ}\cdot\cO^{\wedge+}_{\Spa\,\underline X,x}$ of
$\cO^{\wedge+}_{\Spa\,\underline X,x}$).
\item
The pairs
$(\cO^{\he+}_{\Spa\,\underline X,x},\cO^{\he\,\circ\circ}_{\Spa\,\underline X,x})$
and
$(\cO^{\wedge+}_{\Spa\,\underline X,x},\cO^{\wedge\,\circ\circ}_{\Spa\,\underline X,x})$
are henselian.
\item
Suppose that $x$ is analytic, and let
$$
\fp^\loc\subset\kappa(x^\loc)^+
\qquad
\fp^\he\subset\kappa(x^\he)^+
\qquad
\fp^\wedge\subset\kappa(x^\wedge)^+
$$
be the (unique) prime ideals of height one. Then :
\begin{enumerate}
\item
$\cO^{\loc\,\circ\circ}_{\Spa\,\underline X,x}=\pi^{\loc\,-1}_x(\fp^\loc)
\qquad
\cO^{\he\,\circ\circ}_{\Spa\,\underline X,x}=\pi^{\he\,-1}_x(\fp^\he)
\qquad
\cO^{\wedge\circ\circ}_{\Spa\,\underline X,x}=\pi^{\wedge\,-1}_x(\fp^\wedge)$.
\item
The pairs $(\kappa(x^\he)^+,\fp^\he)$ and
$(\kappa(x^\wedge)^+,\fp^\wedge)$ are henselian.
\item
The local rings $\cO^\he_{\Spa\,\underline X,x}$ and
$\cO^\wedge_{\Spa\,\underline X,x}$ are henselian.
\item
The valuation rings $\kappa(x^\he)^+_{\fp^\he}$ and
$\kappa(x^\wedge)^+_{\fp^\wedge}$ are henselian.
\end{enumerate}
\item
Suppose that $x$ is non-analytic. Then :
\begin{enumerate}
\item
$\cO^{\loc\,\circ\circ}_{\Spa\,\underline X,x}$ (resp.
$\cO^{\he\,\circ\circ}_{\Spa\,\underline X,x}$, resp.
$\cO^{\wedge\,\circ\circ}_{\Spa\,\underline X,x}$) is also an ideal of
$\cO^\loc_{\Spa\,\underline X,x}$ (resp. of $\cO^\he_{\Spa\,\underline X,x}$,
resp. of $\cO^\wedge_{\Spa\,\underline X,x}$), and it is the radical of
the ideal generated by $A^{\circ\circ}$.
\item
The pairs
$(\cO^\he_{\Spa\,\underline X,x},\cO^{\he\,\circ\circ}_{\Spa\,\underline X,x})$
and
$(\cO^\wedge_{\Spa\,\underline X,x},\cO^{\wedge\circ\circ}_{\Spa\,\underline X,x})$
are henselian.
\end{enumerate}
\end{enumerate}
\end{lemma}
\begin{proof} Let $f\in\cO^\wedge_{\Spa\,\underline X,x}$, and pick
any rational subset $R$ in $\Spa\underline X$ with $x\in R$,
and such that $f$ extends to a section
$f_R\in B:=\cO^\wedge_{\Spa\,\underline X}(R)$. Let
$\underline Y^\wedge:=(Y^\wedge,\cT^\wedge_Y,A^{\wedge+}_Y)$ be a
complete and separated quasi-affinoid scheme representing
the subsheaf $h''_R$ of $h''_{\underline X^\wedge}$, so that
$B=\cO_{Y^\wedge}(Y^\wedge)$.

(i): Suppose now that $|f|^\wedge_x\neq 0$; then clearly
$v^\wedge_{x,R}(f_R)\neq 0$. We set
$U:=R_B\bigl(\frac{f_R}{f_R}\bigr)\cap\Spa\,\underline Y^\wedge$.
In view of lemma \ref{lem_invariance-by-loc-hens}(iii) and
corollary \ref{cor_rational-represent}(i), we may then identify
$U$ with an open subset of $\Spa\,X$ containing $x$, and we
pick a rational subset $U'\subset U$ of $\Spa\,\underline X$
with $x\in U'$; by construction the image of $f_R$ under the
restriction map $\rho_{U'}:\cO^\wedge_{\Spa\,\underline X}(R)\to
\cO^\wedge_{\Spa\,\underline X}(U')$ is an invertible element.
This shows that $\Ker\,|\cdot|^\wedge_x$ is the unique maximal
ideal of $\cO^\wedge_{\Spa\,\underline X,x}$, and it also easily
implies that the projection
$\pi^\wedge_x:\cO^\wedge_{\Spa\,\underline X,x}\to\kappa(x^\wedge)$
is surjective (details left to the reader). The same argument
applies to $\cO^\he_{\Spa\,\underline X,x}$ and
$\cO^\loc_{\Spa\,\underline X,x}$, whence the contention.

(ii): Suppose next that $|f|^\wedge_x\leq 1$. Then we let
$U:=R_B\bigl(\frac{f_R}{1}\bigr)\cap\Spa\,\underline Y^\wedge$,
which again we identify naturally with an open subset of
$\Spa\,\underline X$ containing $x$, and we pick a rational
subset $U'\subset U$ of $\Spa\,\underline X$ with $x\in U'$;
then $\rho_{U'}(f_R)\in\cO^{\wedge+}_{\Spa\,\underline X}(U')$. This
shows that $f\in\cO^{\wedge+}_{\Spa\,\underline X,x}$, whence the
assertion concerning the third square diagram. The same
argument applies to the other two diagrams. It follows easily
that the preimage in $\cO^{\wedge+}_{\Spa\,\underline X,x}$ of the
maximal ideal of $\kappa(x^\wedge)^+$ is the unique maximal
ideal of $\cO^{\wedge+}_{\Spa\,\underline X,x}$ (details left to
the reader); especially, the latter is a local ring. Again,
the same applies to the other two rings.

(iii): It follows easily from corollary
\ref{cor_rational-represent}(iv) that the natural map
$A\to\cO^\loc_{\Spa\,\underline X}(U)$ is f-adic, for every
rational subset $U\subset\Spa\,\underline X$. Hence,
$A^{\circ\circ}\cdot\cO^{\loc+}_{\Spa\underline X}(U)$ is an
open ideal of $\cO^{\loc+}_{\Spa\underline X}(U)$ consisting
of topologically nilpotent section, and therefore its
radical equals $\cO^{\loc\,\circ\circ}_{\Spa\underline X}(U)$.
It suffices then to observe that the radical of
$A^{\circ\circ}\cdot\cO^{\loc+}_{\Spa\underline X,x}$ is the
colimit of the filtered system of the radicals of the
ideals
$(A^{\circ\circ}\cdot\cO^{\loc+}_{\Spa\underline X}(U)~|~x\in U)$.
The same applies to the other ideals in (iii).

(iv): By proposition \ref{prop_quasi-affinoid}(iii), both
$(\cO^{\he+}_{\Spa\,\underline X}(U),
\cO^{\he\,\circ\circ}_{\Spa\,\underline X}(U))$ and
$(\cO^{\wedge+}_{\Spa\,\underline X}(U),
\cO^{\wedge\circ\circ}_{\Spa\,\underline X}(U))$ are henselian pairs
for every rational subset $U\subset\Spa\,\underline X$
containing $x$. The assertion follows, after taking colimits
over the filtered family of such rational subsets.

(v.a): First, recall that the existence of a rank one
prime ideal $\fp^\loc\subset\kappa(x^\loc)^+$ is ensured
by lemma \ref{lem_Cont-A}(v); it follows easily that
$\fp^\loc=\kappa(x^\loc)^{\circ\circ}$, and taking into
account remark \ref{rem_Cont-A}(i) we deduce already
that $\pi^\loc_x(\cO^{\loc\,\circ\circ}_{\Spa\,\underline X,x})
\subset\fp^\loc$. Conversely, let
$f\in\cO^\loc_{\Spa\,\underline X,x}\setminus
\cO^{\loc\,\circ\circ}_{\Spa\,\underline X,x}$ be any element;
denote by $\cU$ the set of all rational subsets of
$\Spa\,\underline X$ containing $x$, and for every
$U\in\cU$ set $U':=U\cap R_A\bigl(\frac{1}{f}\bigr)$.
The assumption on $f$ implies that $U'$ is a non-empty
constructible subset of $\Spa\,\underline X$ for every
$U\in\cU$ (corollary \ref{cor_corcor}(ii)), and therefore
$T:=\bigcap_{U\in\cU}U'\neq\emptyset$ (proposition
\ref{prop_lim-qc-sober}(ii.a)). However, $T$ lies in the
set $\Spa\,\underline X(x)$ of all generization of $x$
in $\Spa\,\underline X$ (remark \ref{rem_specialize}(i)),
and must then contain the unique maximal generization of
$x$; the latter means precisely that
$\pi^\loc_x(f)\notin\fp^\loc$, as required.

Next, since $|\cdot|^\loc_x$ and $|\cdot|^\wedge_x$ have the
same value group, the existence of a (unique) height one
prime ideal $\fp^\wedge\subset\kappa(v^\wedge)^+$ follows
from the existence of $\fp^\loc$ and remark
\ref{rem_valuations}(vii); then, since $v^\wedge_{x,R}$ is
continuous, we deduce likewise that
$\pi^\wedge_x(\cO^{\wedge\circ\circ}_{\Spa\,\underline X,x})\subset
\fp^\wedge$ for every rational subset
$R\subset\Spa\,\underline X$. If now
$f\in\cO^\wedge_{\Spa\,\underline X,x}\setminus
\cO^{\wedge\circ\circ}_{\Spa\,\underline X,x}$, pick $R$ as in the
foregoing and $\underline Y$ representing $h''_R$, so that
$x\in R$, and $f$ extends to a section
$f_R\in B:=\cO^\wedge_{\Spa\,\underline Y}(Y)$; let $\cU$ be the
set of all rational subsets of $\Spa\,\underline Y$
containing $x$, and for every $U\in\cU$ set
$U':=U\cap R_B\bigl(\frac{1}{f_R}\bigr)$. Due to corollary
\ref{cor_rational-represent}(i) and lemma
\ref{lem_invariance-by-loc-hens}(iii), the image of
$T:=\bigcap_{U\in\cU}U'$ in $\Spa\,\underline X$ lies in
$\Spa\,\underline X(x)$, and is not empty, and arguing
as in the foregoing we deduce again that
$\pi^\wedge_x(f)\notin\fp^\wedge$. The same argument applies
to $\fp^\he$.

(v.b): Let us remark more generally :

\begin{claim}\label{cl_transit-henselian}
Let $A$ be any ring, $J\subset I\subset A$ two ideals.
Then the pair $(A,I)$ is henselian if and only if the
same holds for both $(A/J,I/J)$ and $(A,J)$.
\end{claim}
\begin{pfclaim} From \cite[Rem.5.1.10(iv,v)]{Ga-Ra} we
see that if $(A,I)$ is henselian, the same holds for the
pairs $(A,J)$ and $(A/J,I/J)$. Conversely, let $A\to B$
be a finite ring homomorphism; if $(A,J)$ is henselian,
the projection $B\to B/JB$ restricts to a bijection
$\mathrm{Idemp}\,(B)\isom\mathrm{Idemp}\,(B/JB)$ on the
respective subsets of idempotent elements. Then, if also
$(A/J,I/J)$ is henselian, the induced map
$\mathrm{Idemp}\,(B/JB)\to\mathrm{Idemp}\,(B/IB)$ is bijective
as well. Summing up, we conclude that $(A,I)$ is henselian.
\end{pfclaim}

In view of (iv) and claim \ref{cl_transit-henselian},
we are easily reduced to showing that the projections
$\cO^{\he+}_{\Spa\,\underline X,x}\to\kappa(x^\he)^+$ and
$\cO^{\wedge+}_{\Spa\,\underline X,x}\to\kappa(x^\wedge)^+$ are
surjective. The latter follows from (ii) and from the
surjectivity of $\pi^\he_x$ and $\pi^\wedge_x$, which has
already been remarked in the proof of (i).

(v.c): Let $\fm^\he:=\Ker\,\pi^\he_x$; the proof of (i) shows
that $\fm^\he$ is the maximal ideal of
$\cO^\he_{\Spa\,\underline X,x}$, and (ii) implies that
$\fm^\he$ is also an ideal in $\cO^{\he+}_{\Spa\,\underline X,x}$;
moreover, in light of (iv) and \cite[Rem.5.1.10(iv)]{Ga-Ra}
the pair $(\cO^{\he+}_{\Spa\,\underline X,x},\fm^\he)$ is henselian.
Then the assertion follows from \cite[Rem.5.1.10(ii)]{Ga-Ra}.
The same argument applies to $\cO^\wedge_{\Spa\,\underline X,x}$.

(v.d): It is easily seen that the localization map
$\kappa(x^\he)^+\to\kappa(x^\he)^+_{\fp^\he}$ restricts to
a bijection from $\fp^\he$ to the maximal ideal of
$\kappa(x^\he)^+_{\fp^\he}$; then the assertion follows
from \cite[Rem.5.1.10(ii)]{Ga-Ra}. The same argument
applies to $\kappa(x^\wedge)^+_{\fp^\wedge}$.

(vi.a): By (iii) we know already that
$\cO^{\wedge\circ\circ}_{\Spa\,\underline X,x}$ lies in the radical
of $A^{\circ\circ}\cdot\cO^\wedge_{\Spa\,\underline X,x}$; hence,
suppose that $f^n\in A^{\circ\circ}\cdot\cO^\wedge_{\Spa\,\underline X,x}$
for some $n\in\N$, and say that $f^n=\sum_{i=1}^ka_ig_i$ for
some $a_1,\dots,a_k\in A^{\circ\circ}$ and
$g_1,\dots,g_k\in\cO^\wedge_{\Spa\,\underline X,x}$. Notice that
since $x$ is non-analytic, we have $\pi^\wedge_x(a_i)=0$
for $i=1,\dots,k$, whence $\pi^\wedge_x(f^n)=0$, and therefore
$\pi^\wedge_x(f^ng_i)=0$ as well for $i=1,\dots,k$; then (ii)
implies that $f^ng_i\in\cO^{\wedge+}_{\Spa\,\underline X,x}$, so
$f^{2n}=\sum_{i=1}^ka_i\cdot(f^ng_i)\in
A^{\circ\circ}\cdot\cO^{\wedge+}_{\Spa\,\underline X,x}$. The same
argument applies to $\cO^{\loc\,\circ\circ}_{\Spa\,\underline X,x}$
and $\cO^{\he,\circ\circ}_{\Spa\,\underline X,x}$.

(vi.b) follows directly from (vi.a), (iv) and
\cite[Rem.5.1.10(ii)]{Ga-Ra}.
\end{proof}

\sset\subsubsection{}\label{subsec_ZR}
Let $\underline A:=(A,A^+,U)$ be any quasi-affinoid ring,
and $J\subset A$ a finitely generated ideal such that
$U=\Spec\,A\setminus\Spec\,A/J$. The {\em Zariski-Riemann
spectrum} of $\underline A$ is the subset of $\Spv\,A$
$$
\ZR(\underline A):=\Spa(\underline A)\cap\Spv(A,J)
$$
(notation of definition \ref{def_Spv-I}(i)) that we endow
with the topology induced by the inclusion into $\Spv\,A$.
A {\em rational subset\/} of $\ZR(\underline A)$ is a subset
of the form $R\cap\ZR(\underline A)$, where $R$ is a rational
subset of $\Spv(A,J)$ (see definition \ref{def_Spv-I}(ii)).
Notice that neither $\ZR(\underline A)$ nor the class of
its rational subsets depend on the choice of $J$ (lemma
\ref{lem_was-rem-Spv-A-I}(i)).

\begin{proposition}\label{prop_ZR}
With the notation of \eqref{subsec_ZR}, the following holds :
\begin{enumerate}
\item
$\ZR(\underline A)=\{v\in\Spa(\underline A)~|~
\text{$v$ has no proper primary specializations in
$\Spa(\underline A)$}\}$.
\item
$\ZR(\underline A)$ is a pro-constructible subset of\/ $\Spv(A,J)$.
In particular, it is a spectral topological space.
\item
The rational subsets of\/ $\ZR(\underline A)$ are constructible
in $\ZR(\underline A)$, and form a basis of the topology of\/
$\ZR(\underline A)$ that is closed under finite intersections.
\item
The unit of adjunction
$\underline A\to\sGamma\circ\sSpec\,\underline A$ induces
a homeomorphism
$$
\ZR(\sGamma\circ\sSpec\,\underline A)\isom\ZR(\underline A).
$$
Especially, $\ZR(\underline A)$ depends only on the
quasi-affinoid scheme $\sSpec\,\underline A$.
\end{enumerate}
\end{proposition}
\begin{proof}(i): Recall that $\Cont(A)$ is closed under
specializations in $\Spv\,A$ (remark \ref{rem_Cont-A}(iv));
it follows easily that $\Spa(A,A^+)$ is closed under primary
specializations.
Recall as well that a proper specialization $(v,w)$ in
$\Spv\,A$ is $J$-admissible if and only if $w$ is a primary
specialization of $v$ and the supports of $v$ and $w$ lie
in $U$ (see \eqref{subsec_I-admissible}); hence, if
$v\in\Spa\,\underline A$, the proper specialization
$(v,w)$ is $J$-admissible if and only if $w$ is a primary
specialisation of $v$ that lies in $\Spa\,\underline A$.
Then the assertion follows from lemma \ref{lem_cGamma-I}(iii).

(ii): Let $A_0\subset A$ be a subring of definition, and
$I\subset A_0$ a finitely generated ideal of adic definition;
recall that $\Spa(\underline A)$ is a pro-constructible
subset of $\Cont(A)$ (proposition \ref{prop_f-adic-Spectra}(i))
and the inclusion map $\Cont(A)\to\Spv(A,IA)$ is a closed
immersion (theorem \ref{th_Cont-spectral}(ii)) so
$\Spv(\underline A)$ is a pro-constructible subset of
$\Spv(A,IA)$ (corollary \ref{cor-pro-constr}(i)).
Moreover, the radical of $J$ contains $IA$, hence the
inclusion map $\Spv(A,J)\to\Spv(A,IA)$ admits a spectral
retraction (lemma \ref{lem_was-rem-Spv-A-I}(i))
\set\begin{equation}\label{eq_restrict-retract}
r:\Spv(A,IA)\to\Spv(A,J).
\end{equation}
Namely, for every $v\in\Spv(A,IA)$, the valuation $r(v)$
is the unique $J$-admissible specialization of $v$ that
does not admit further proper $J$-admissible specializations.
We claim that
$$
\ZR(\underline A)=r(\Spa\,\underline A).
$$
Indeed, clearly $\ZR(\underline A)=r(\ZR(\underline A))\subset
r(\Spa\,\underline A)$. The converse inclusion follows from
the proof of (i). Then, the assertion follows from corollary
\ref{cor-pro-constr}(i).

(iii) follows immediately from (ii) and theorem \ref{th_Spv-I}(ii).

(iv) follows directly from (i) and lemmata
\ref{lem_invariance-by-loc-hens}(ii) and
\ref{lem_image-of-special}(i).
\end{proof}

\begin{remark} In \cite{Tem}, Temkin associates with every
separated morphism $f:Y\to X$ of quasi-compact and quasi-separated
schemes, a spectral space that he calls the Riemann-Zariski
space of $f$, and denotes $\mathrm{RZ}_Y(X)$, together with a
continuous map $\mathrm{RZ}_Y(X)\to X$. His construction generalizes
earlier work by Zariski, and is related as follows to our
$\ZR(\underline A)$. Let $U\to\Spec\,A^+$ be the composition
of the open immersion $U\to\Spec\,A$ with the natural morphism
$\Spec\,A\to\Spec\,A^+$; then we may consider the induced map
$\phi:\mathrm{RZ}_U(\Spec\,A^+)\to\Spec\,A^+$, and there is a
natural homeomorphism
$\phi^{-1}(\Spec\,A^+/A^{\circ\circ})\isom\ZR(\underline A)$.
\end{remark}

\sset\subsubsection{}\label{subsec_standard-cov}
Keep the notation of \eqref{subsec_ZR}, and for some integer
$n\geq -1$, let $R_\bullet:=(R_i~|~i=0,\dots,n)$ be a finite
system of open subsets of $\ZR(\underline A)$ such that
$\ZR(\underline A)=\bigcup^n_{i=0}R_i$ (the case $n=-1$
occurs for the empty covering of the empty Zariski-Riemann
spectrum). Then we say that $R_\bullet$ is a
{\em standard covering\/} of $\ZR(\underline A)$ if there
exists a finite system $f_\bullet:=(f_0,\dots,f_n)$ of
elements of $A$ such that :
\begin{itemize}
\item
$R_i=R_A\bigl(\frac{f_0}{f_i},\dots,\frac{f_n}{f_i}\bigr)
\cap\ZR(\underline A)$ for $i=0,\dots,n$
\item
$J$ is contained in the radical of the ideal of $A$ generated
by the system $f_\bullet$.
\end{itemize}
Especially, each $R_i$ is a rational subset of $\ZR(\underline A)$.
By construction, the retraction $r$ of \eqref{eq_restrict-retract}
restricts to a surjective continuous map
$$
s:\Spa\,\underline A\to\ZR(\underline A)
$$
and if $f_\bullet$ fulfills the foregoing conditions, lemma
\ref{lem_characterize-I-closed}(i) says that
$$
s^{-1}R_i=R_A\Bigl(\frac{f_0}{f_i},\dots,\frac{f_n}{f_i}\Bigr)
\cap\Spa\,\underline A
\qquad
\text{for every $i=0,\dots,n$}.
$$
Especially, the system
$(R_A\bigl(\frac{f_0}{f_i},\dots,\frac{f_n}{f_i}\bigr)\cap
\Spa\,\underline A~|~i=0,\dots,n)$ is a covering of
$\Spa\,\underline A$. A covering of $\Spa\,\underline A$
of this type shall also be called a {\em standard covering}.

\begin{lemma}\label{lem_standard-coverings}
If $A$ is topologically local, every open covering of\/
$\ZR(\underline A)$ or $\Spa\,\underline A$ can be refined
by a standard covering.
\end{lemma}
\begin{proof} Let $(U_\lambda~|~\lambda\in\Lambda)$ be an
open covering of $\ZR(\underline A)$; we need to find a
standard covering $(R_i~|~i=0,\dots,n)$ such that, for every
$i=0,\dots,n$ there exists $\lambda\in\Lambda$ with
$R_i\subset U_\lambda$. To this aim, in light of proposition
\ref{prop_ZR}(i,ii) we may assume that $\Lambda$ is a finite
set, and each $U_\lambda$ is rational, say
$$
U_\lambda=R_A\Bigl(\frac{g_{\lambda,1}}{g_{\lambda,0}},\dots,
\frac{g_{\lambda,n_\lambda}}{g_{\lambda,0}}\Bigr)\cap\ZR(\underline A)
$$
for a finite system
$g_{\lambda,\bullet}:=(g_{\lambda,j}~|~j=0,\dots,n_\lambda)$
of elements of $A$ that generates an ideal whose radical
contains $J$. Define $S$ as the set of all the sequences
of integers $j_\bullet:=(j_\lambda~|~\lambda\in\Lambda)$
such that $0\leq j_\lambda\leq n_\lambda$ for every
$\lambda\in\Lambda$. Let also $T\subset S$ be the subset
of all sequences $j_\bullet$ such that $j_\lambda=0$ for at
least one $\lambda\in\Lambda$. Set
$g_{j_\bullet}:=\prod_{\lambda\in\Lambda}g_{\lambda,j_\lambda}$ for
every $j_\bullet\in S$. Notice that :
$$
R_{k_\bullet}:=
R_A\Bigl(\frac{g_{j_\bullet}}{g_{k_\bullet}}~|~j_\bullet\in S\Bigr)=
\bigcap_{\lambda\in\Lambda}R_A
\Bigl(\frac{g_{\lambda,1}}{g_{\lambda,k_\lambda}},\dots,
\frac{g_{\lambda,n_\lambda}}{g_{\lambda,k_\lambda}}\Bigr)
\qquad
\text{for every $k_\bullet\in S$}.
$$
Moreover, each $g_{\lambda,\bullet}$ defines a standard open
covering of $\ZR(\underline A)$. It follows easily that
\set\begin{equation}\label{eq_from-FP}
\ZR(\underline A)=\bigcup_{k_\bullet\in T}
(R_{k_\bullet}\cap\ZR(\underline A))
\end{equation}
and furthermore, for every $k_\bullet\in T$, the subset
$R_{k_\bullet}\cap\ZR(\underline A)$ is contained in some $U_\lambda$.

\begin{claim}\label{cl_from-FP}
$R_{k_\bullet}\cap\ZR(\underline A)=
R_A\bigl(\frac{g_{j_\bullet}}{g_{k_\bullet}}~|~j_\bullet\in T\bigr)
\cap\ZR(\underline A)$ for every $k_\bullet\in T$.
\end{claim}
\begin{pfclaim} Indeed, the inclusion
$R_{k_\bullet}\subset R'_{k_\bullet}:=
R_A\bigl(\frac{g_{j_\bullet}}{g_{k_\bullet}}~|~j_\bullet\in T\bigr)$
is obvious. Conversely, let $v\in R'_{k_\bullet}\cap\ZR(A)$;
we see from \eqref{eq_from-FP} that for every
$j_\bullet\in S$ there exists $i_\bullet\in T$ such that
$v(g_{j_\bullet})\leq v(i_\bullet)\neq 0$, whence
$v(g_{j_\bullet})\leq v(k_{\bullet})\neq 0$, so $v\in R_{k_\bullet}$.
\end{pfclaim}

In light of \eqref{eq_from-FP} and claim \ref{cl_from-FP},
to conclude it now suffices to show :

\begin{claim} $J$ is contained in the radical of the ideal $J'$
of $A$ generated by $(g_{k_\bullet}~|~k_\bullet\in T)$.
\end{claim}
\begin{pfclaim} First, we check that $J'$ is an open
ideal, using the criterion of lemma \ref{lem_criterion-opennes}.
Indeed, by claim \ref{cl_rank-one-vals}, every continuous
analytic rank one valuation of $A$ lies in $\Spa\,\underline A$,
and from \eqref{eq_from-FP} and claim \ref{cl_from-FP} we see
that for every such valuation $v$ there exists $k_\bullet\in T$
such that $v(g_{k_\bullet})\neq 0$, whence the assertion.
By lemma \ref{lem_deja-vu}(v), it follows that
$\Spec A/J'\subset X^{\circ\circ}_A:=\Spec\,A/A^{\circ\circ}A$.
Hence, let $\fp\in X^{\circ\circ}_A\setminus\Spec A/J'$; it
suffices to show that $J'\not\subset\fp$. To this aim, denote
by $v_\fp$ the trivial valuation of $A$ supported at $\fp$
(see remark \ref{rem_semi-norm}(vii)); since $\fp$ is an open
ideal (lemma \ref{lem_deja-vu}(i)), the valuation $v_\fp$ is
continuous (remark \ref{rem_Cont-A}(vi)), and then obviously
$v\in\Spa\,\underline A$. As in the foregoing, we deduce that
$v_\fp(g_{k_\bullet})\neq 0$ for some $k_\bullet\in T$, {\em i.e.}
$g_{k_\bullet}\notin\fp$, whence the claim.
\end{pfclaim}

Next, consider any open covering $(U'_\lambda~|~\lambda\in\Lambda)$
of $\Spa\,\underline A$, and set
$U_\lambda:=U'_\lambda\cap\ZR(\underline A)$ for every
$\lambda\in\Lambda$; clearly $(U_\lambda~|~\lambda\in\Lambda)$
is an open covering of $\ZR(\underline A)$, which can be
refined by a standard covering $(R_0,\dots,R_n)$; but then the
system $(s^{-1}R_0,\dots,s^{-1}R_n)$ is a standard covering of
$\Spa\,\underline A$, and notice that if $R_i\subset U_\lambda$,
we get
$$
s^{-1}R_i\subset s^{-1}U_\lambda\subset U'_\lambda
$$
where the last inclusion holds, since $U'_\lambda$ is open in
$\Spa\,\underline A$ and $s^{-1}U_\lambda$ is the set of all primary
generizations of the elements of $U_\lambda$ in $\Spa\,\underline A$.
Thus, $(s^{-1}R_0,\dots,s^{-1}R_n)$ refines the covering
$(U'_\lambda~|~\lambda\in\Lambda)$, as required.
\end{proof}

\begin{remark}\label{rem_extract-proj}
The following constructions shall be exploited in the proofs
of both theorem \ref{th_henselians-are-sheaves} and theorem
\ref{th_compl-an-noeth-are-sheaves}.

(i)\ \
Let $\underline X$ be a topologically local quasi-affinoid
scheme, set $(A,A^+,X):=\sGamma(\underline X)$, and choose a
subring of definition $A_0\subset A$ and a finitely generated
ideal $I_0\subset A_0$ of adic definition. Consider a sequence
$f_\bullet:=(f_0,\dots,f_n)$ of elements of $A$ that generates
an ideal $J\subset A$ such that $X\cap\Spec\,A/J=\emptyset$,
and let $R_\bullet:=(R_0,\dots,R_n)$ be the standard covering
of $\Spa\,\underline X$ associated with $f_\bullet$. For every
subset $\Lambda\subset\{0,\dots,n\}$, let
$A_\Lambda:=A[f_i^{-1}~|~i\in\Lambda]$ and
$X_\Lambda:=X\cap\Spec\,A_\Lambda$; by inspecting the definitions
we see that
$$
R_\Lambda:=\bigcap_{i\in\Lambda}R_i=\Spa\,\underline A{}_\Lambda
\qquad\text{with}\qquad
\underline A{}_\Lambda:=(A_\Lambda,A^+_\Lambda,X_\Lambda)
$$
where $A^+_\Lambda$ is the integral closure of
$A^+[f_k/f_i~|~(k,i)\in\{0,\dots,n\}\times\Lambda]$ in
$A_\Lambda$; moreover,
$A_{0,\Lambda}:=A_0[f_k/f_i~|~(k,i)\in\{0,\dots,n\}\times\Lambda]$
is a subring of definition of $A_\Lambda$, and the natural ring
homomorphism $A_0\to A_{0,\Lambda}$ is adic (here we let
$R_\emptyset:=\Spa\,\underline X$). Furthermore, since
$X$ is quasi-compact, we have $\cO_{\!X}(X_\Lambda)=A_\Lambda$
(\cite[Ch.I, Prop.9.2.1]{EGAI}), and therefore
$$
\underline A{}_\Lambda=\sGamma\circ\sSpec\,(\underline A{}_\Lambda)
\qquad
\text{for every $\Lambda\subset\{0,\dots,n\}$}.
$$

(ii)\ \
Denote by $(A_{0,\Lambda}^\he,I^\he_{0,\Lambda})$ the henselization
of the pair $(A_{0,\Lambda},I_0A_{0,\Lambda})$, and endow
$A_{0,\Lambda}^\he$ with the $I^\he_{0,\Lambda}$-adic topology; we
deduce a unique isomorphism of topological $A$-algebras :
$$
\omega^\he_\Lambda:
A_\Lambda^\he:=A_\Lambda\otimes_{A_{0,\Lambda}}A_{0,\Lambda}^\he\isom
\cO^\he_{\Spa\,\underline X}(R_\Lambda)
\qquad
\text{for every $\Lambda\subset\{0,\dots,n\}$}
$$
where the topology on $A^\he_\Lambda$ is the unique one such
that the natural map $A_{0,\Lambda}^\he\to A_\Lambda^\he$ is open
(proposition \ref{prop_f-adic-push-out}(ii)). Likewise, let
$A_\Lambda^\wedge$ be the separated completion of $A_\Lambda$;
we deduce a natural continuous homomorphism of $A$-algebras
$$
\omega^\wedge_\Lambda:A_\Lambda^\wedge\to
\cO^\wedge_{\Spa\,\underline X}(R_\Lambda) 
\qquad
\text{for every $\Lambda\subset\{0,\dots,n\}$}.
$$
Notice also that for every $\Lambda'\subset\Lambda$ there
exist unique continuous homomorphisms of $A$-algebras and
respectively $A_{0,\Lambda}$-algebras
$$
\mu_{\Lambda',\Lambda}:A_{\Lambda'}\to A_\Lambda
\qquad
\mu^\he_{0,\Lambda',\Lambda}:A_{0,\Lambda'}^\he\to A_{0,\Lambda}^\he
$$
and set $\mu^\he_{\Lambda',\Lambda}:=
\mu_{\Lambda',\Lambda}\otimes_A\mu^\he_{0,\Lambda',\Lambda}$. Let also
$\mu^\wedge_{\Lambda',\Lambda}:A^\wedge_{\Lambda'}\to A^\wedge_\Lambda$
be the completion of $\mu_{\Lambda',\Lambda}$, and
$$
\rho^\he_{\Lambda',\Lambda}:\cO^\he_{\Spa\,\underline X}(R_{\Lambda'})
\to\cO^\he_{\Spa\,\underline X}(R_\Lambda)
\qquad
\rho^\wedge_{\Lambda',\Lambda}:\cO^\wedge_{\Spa\,\underline X}(R_{\Lambda'})
\to\cO^\wedge_{\Spa\,\underline X}(R_\Lambda)
$$
the restriction homomorphisms associated with the inclusion
$R_\Lambda\subset R_{\Lambda'}$.

(iii)\ \
We claim that the diagrams of $A$-algebras:
$$
\xymatrix{
A^\he_{\Lambda'} \ar[r]^-{\omega^\he_{\Lambda'}}
\ar[d]_{\mu^\he_{\Lambda',\Lambda}} &
\cO^\he_{\Spa\,\underline X}(R_{\Lambda'})
\ar[d]^{\rho^\he_{\Lambda',\Lambda}} &
A^\wedge_{\Lambda'} \ar[r]^-{\omega^\wedge_{\Lambda'}}
\ar[d]_{\mu^\wedge_{\Lambda',\Lambda}} &
\cO^\wedge_{\Spa\,\underline X}(R_{\Lambda'})
\ar[d]^{\rho^\wedge_{\Lambda',\Lambda}} \\
A^\he_\Lambda \ar[r]^-{\omega^\he_\Lambda} &
\cO^\he_{\Spa\,\underline X}(R_\Lambda) &
A^\wedge_\Lambda \ar[r]^-{\omega^\wedge_\Lambda} &
\cO^\wedge_{\Spa\,\underline X}(R_\Lambda)
}$$
commute for every $\Lambda'\subset\Lambda\subset\{0,\dots,n\}$.
Indeed, a simple inspection shows that $\mu^\he_{\Lambda',\Lambda}$
is continuous. On the other hand, all the arrows in the diagram
are $A$-algebra homomorphisms, therefore also $A_{\Lambda'}$-algebra
homomorphisms, and by the universal property of the topological
henselization there is a unique continuous $A_{\Lambda'}$-algebra
homomorphism $A_{\Lambda'}^\he\to\cO^\he_{\Spa\,\underline X}(R_\Lambda)$,
whence the contention for the left diagram. For the right diagram
one argues similarly, using the universal property of the
separated completions.

(iv)\ \
Endow $A[Y]$ and $B:=A_0[T_0,\dots,T_n]$ with the standard
$\N$-gradings such that $\gr_kA[Y]=AY^k$ and $\gr_kB$ is the
free $A$-module generated by the monomials
$T_0^{\nu_0}\cdots T_n^{\nu_n}$ of total degree
$\nu_0+\cdots+\nu_n=k$, for every $k\in\N$. Set as well
$$
S:=\Spec\,A_0
\qquad
\P^n_S:=\Proj\,B
$$
and recall that $\P^n_S$ admits the following finite affine
open covering (see \eqref{subsec_projective-spectra}).
For every $i=0,\dots,n$ we have the subring
$B_i:=A_0[T_0/T_i,\dots,T_n/T_i]$ of $B[1/T_i]$, and
$\Omega_i:=\Spec\,B_i$ is naturally identified with the
open subscheme of $\P^n_S$ consisting of all homogeneous
prime ideals of $B$ that do not contain $T_i$; then
$$
\P^n_S=\Omega_0\cup\cdots\cup\Omega_n.
$$
Next, according to \eqref{subsec_in-the-situat}, the
homomorphism of $\N$-graded $A_0$-algebras
$$
h:B\to A[Y]
\qquad
T_j\mapsto f_jY
\qquad
\text{for every $j=0,\dots,n$}
$$
induces a morphism of $A_0$-schemes
$$
\phi:\Spec\,A\setminus\Spec\,A/J\to\P^n_S.
$$
Explicitly, $\phi^{-1}\Omega_i=\Spec\,A_{\{i\}}$ for every
$i=0,\dots,n$, and the restriction $\phi^{-1}\Omega_i\to\Omega_i$
of $\phi$ corresponds to the homomorphism of $A_0$-algebras
$$
h_i:B_i\to A_{\{i\}}
\qquad
T_j/T_i\mapsto f_j/f_i
\qquad
\text{for every $j=0,\dots,n$}.
$$

(v)\ \
Now, quite generally, let $\beta:\sX\to\sY$ be any quasi-compact
and quasi-separated morphism of schemes. Then $\beta_*\cO_\sX$
is a quasi-coherent $\cO_\sY$-module (\cite[Ch.I, Prop.9.2.1]{EGAI}),
hence $\cI:=\Ker\,\beta^\sharp:\cO_\sY\to\psi_*\cO_\sX$ is a
quasi-coherent sheaf of ideals of $\cO_\sY$, and therefore
$\Spec\,\cO_\sY/\cI$ is a well defined closed subscheme of
$\sY$, called the {\em schematic image of\/ $\psi$}.
From the definition, it is clear that the construction of the
schematic image is local on $\sY$ : namely, if $\sY'\subset\sY$
is any open subscheme, and $\sX':=\beta^{-1}\sY'$, then the
schematic image of $\beta_{|\sX'}:\sX'\to\sY'$ is the intersection
of $\sY'$ with the schematic image of $\beta$.

Thus, if $V$ is the schematic image of our $\phi$, we get a
commutative diagram of schemes
$$
\xymatrix{ X \ar[rr]^-\psi \ar[d] & & V \ar[d]^\iota \\
\Spec\,A\setminus\Spec\,A/J \ar[rr]^-\phi \ar[rru] & & \P^n_S
}$$
whose left vertical arrow is an open immersion, and with
$\iota$ a closed immersion. Lastly, set $V_\emptyset:=V$
and for every subset $\Lambda\subset\{0,\dots,n\}$ let
$V_\Lambda:=\iota^{-1}\Bigl(\bigcap_{i\in\Lambda}\Omega_i\Bigr)$;
then $V_\Lambda$ is naturally identified with
$\Spec\,A_{0,\Lambda}$, for every non-empty $\Lambda$.

(vi)\ \
Furthermore, let $U$ be the analytic locus of $S$; recall
that the morphism $\Spec\,A\to S$ (associated with the
inclusion map $A_0\to A$) identifies $U$ with the analytic
locus of $\Spec\,A$ (lemma \ref{lem_deja-vu}(iii)), and
$U$ is an open subset of $X$, under this identification.
Then $\P^n_U:=\P^n_S\times_SU$ is an open subset of $\P^n_S$,
and we claim that $\psi$ restricts to an isomorphism of schemes
$$
\psi_{|U}:U\isom V\cap\P^n_U.
$$
Indeed, $\psi_{|U}$ is a closed immersion, by
\cite[Ch.I, Cor.5.4.6]{EGAI}. On the other hand, $\psi_{|U}$
is also the schematic closure of the restriction
$U\to\P^n_U$ of $\phi$, whence the claim.
\end{remark}

\begin{theorem}\label{th_henselians-are-sheaves}
For any topologically local quasi-affinoid scheme
$\underline X:=(X,\cT_X,A^+)$, the presheaves
$\cO^\he_{\Spa\,\underline X}$ and $\cO^{\he+}_{\Spa\,\underline X}$
are sheaves of rings on $(\cQ(\underline X),J_\cQ)$.
\end{theorem}
\begin{proof} Notice first that -- by virtue of corollary
\ref{cor_corcor}(i) -- if $\cO^\he_{\Spa\,\underline X}$ is a sheaf,
the same holds for $\cO^{\he+}_{\Spa\,\underline X}$, so it
suffices to check the assertion for $\cO^\he_{\Spa\,\underline X}$.
To this aim, for every $U\in\Ob(\cQ(\underline X))$ and every
$\cS\in J_\cQ(U)$, let $h_\cS$ be the corresponding sub-presheaf
of the presheaf $h_U$ on $\cQ(\underline X)$ represented by $U$.
According to claim \ref{cl_plus-construction}(ii,iii), it
suffices to show that the natural map
$$
r_U:\cO^\he_{\Spa\,\underline X}(U)\to\colim_{\cS\in J_\cQ(U)}
\Hom_{\cQ(\underline X)^\wedge}(h_\cS,\cO^\he_{\Spa\,\underline X})
$$
is a bijection for every such $U$. However, let
$\phi_{\underline T/\underline X}:\underline T\to\underline X$
be any morphism of topologically local quasi-affinoid
schemes representing the inclusion $h_U\subset h_{\underline X}$
(see \eqref{subsec_quasi-affoid-subsets}); by remark
\ref{rem_yoneda-rationals}(v) it follows that
$\phi_{\underline T/\underline X}$ induces an equivalence
of categories
$$
\cQ(\underline T)\isom\cQ(\underline X)/U
\qquad
V\mapsto\Spa\,\phi_{\underline T/\underline X}(V)\subset U
$$
so we may replace $\underline X$ by $\underline T$, and
reduce to checking that $r_{\Spa\,\underline X}$ is an isomorphism.
Next, lemma \ref{lem_standard-coverings} implies that
the sieves of $\cQ(\underline X)$ generated by the
standard coverings of $\Spa\,\underline X$ form a final
subset of $J_\cQ(\Spa\,\underline X)$. Hence, let $f_\bullet$
and $R_\bullet$ be as in remark \ref{rem_extract-proj}(i);
taking into account the discussion of
\eqref{subsec_interpret-descent}, we are then reduced to
checking that the natural map
$$
\cO^\he_{\Spa\,\underline X}(\Spa\,\underline X)\to
\Equal\Bigl(\prod_{i=0}^n\xymatrix{\cO^\he_{\Spa\,\underline X}(R_i)
\ar@<-.5ex>[r] \ar@<.5ex>[r] &}
\prod_{i,j=0}^n\cO^\he_{\Spa\,\underline X}(R_{\{i,j\}})\Bigr)
$$
is a bijection. Now, say that $(A,A^+,X)=\sGamma(\underline X)$.
Recall that the henselization morphism
$$
\lambda:\underline X^\he:=(X^\he,\cT_X^\he,A^{\he+})\to\underline X
$$
is f-adic and induces a homeomorphism
$\Spa\,\lambda:\Spa\,\underline X^\he\isom\Spa\,\underline X$,
as well as an isomorphism of schemes
$X^\he\isom X\times_AA^\he$; it follows easily that the system
$((\Spa\,\lambda)^{-1}R_i~|~i=0,\dots,n)$ is a rational covering
of $\Spa\,\underline X^\he$, and moreover there is a unique
isomorphism of presheaves of topological $A$-algebras
$$
\cO^\he_{\Spa\,\underline X^\he}\isom
(\Spa\,\lambda)^{-1}\cO^\he_{\Spa\,\underline X}.
$$
We may then replace $\underline X$ by $\underline X^\he$, and
assume from start that $\underline X$ is topologically henselian.

In this situation, let $A_0\subset A$ be a subring of definition,
and $I_0\subset A_0$ an ideal of adic definition; consider the
system of quasi-affinoid rings
$(\underline A{}_\Lambda~|~\Lambda\subset\{0,\dots,n\})$,
their rings of definition $A_{0,\Lambda}\subset A_\Lambda$,
and the systems of their topological henselizations
$(A^\he_\Lambda~|~\Lambda\subset\{0,\dots,n\})$ and
$(A^\he_{0,\Lambda}~|~\Lambda\subset\{0,\dots,n\})$
constructed in remark \ref{rem_extract-proj}(i).
Define also $S$, $V$, $\psi:X\to V$ and $V_\Lambda$
for every subset $\Lambda\subset\{0,\dots,n\}$ as in
remark \ref{rem_extract-proj}(iv,v); set
$S':=\Spec\,A_0/I_0$, $V'_\Lambda:=V_\Lambda\times_SS'$
for every such $\Lambda$, and consider the cartesian
diagram of schemes
$$
\xymatrix{ V' \ar[r]^-{\tau_V} \ar[d]_{\pi'} & V \ar[d]^\pi \\
S' \ar[r]^-{\tau_S} & S
}$$
where $\pi$ is the restriction to $V$ of the projection
$\P^n_S\to S$ and $\tau_S:S'\to S$ is the closed immersion.
For every scheme $\sX$, let $\sX_\et$ be the \'etale
site of $\sX$, and $\cO_{\sX_\et}$ the structure sheaf of rings
on $\sX_\et$. We define the sheaf on $V'_\et$
$$
\cF:=\tau^*_V\psi_*\cO_{\!X_\et}.
$$

\begin{claim}\label{cl_I-am-convinced}
(i)\ \ 
For every $\Lambda\subset\{0,\dots,n\}$ there
exists a commutative diagram of rings
$$
\cD_\Lambda
\quad : \quad
{\diagram A_{0,\Lambda}^\he \ar[r]^-{\omega_{0,\Lambda}} \ar[d] &
\Gamma(V'_\Lambda,\tau^*_V\cO_{V,\et}) \ar[d] \\
A^\he_\Lambda \ar[r]^-{\omega_\Lambda} & \Gamma(V'_\Lambda,\cF)
\enddiagram}$$
whose left vertical arrow is the inclusion map, and whose
right vertical arrow is deduced from the natural morphism
$\cO_{V,\et}\to\psi_*\cO_{\!X,\et}$. Moreover, $\omega_{0,\Lambda}$
is an isomorphism of $A_0$-algebras for every
$\Lambda\neq\emptyset$, and $\omega_\Lambda$ is an isomorphism
of $A$-algebras for every $\Lambda$.
\begin{enumerate}
\addenu
\item
For every $\Lambda'\subset\Lambda$, let
$\rho_{\Lambda',\Lambda}:\Gamma(V'_{\Lambda'},\cF)\to
\Gamma(V'_\Lambda,\cF)$ be the restriction homomorphism
associated with the inclusion $V'_\Lambda\subset V'_{\Lambda'}$.
Then we have a commutative diagram
$$
\xymatrix{ A^\he_{\Lambda'} \ar[r]^-{\omega_{\Lambda'}}
\ar[d]_{\mu^\he_{\Lambda',\Lambda}} & \Gamma(V'_{\Lambda'},\cF)
\ar[d]^{\rho_{\Lambda',\Lambda}} \\
A^\he_\Lambda \ar[r]^-{\omega_\Lambda} &
\Gamma(V'_\Lambda,\cF).
}$$
\end{enumerate}
\end{claim}
\begin{pfclaim}(i): We consider first the case where
$\Lambda=\emptyset$, and notice that $V_\emptyset=V$,
$A^\he_{0,\emptyset}=A_0$ and $A^\he_\emptyset=A$. Under the
current assumptions, the pair $(A_0,I_0)$ is henselian;
by \cite[Exp.XII, Th.5.1(i) and Prop.6.5(i)]{SGA4-3} we
then get for every sheaf $\cG$ on $V_\et$ natural isomorphisms
$$
\Gamma(V,\cG)\isom\Gamma(S,\pi_*\cG)\isom
\Gamma(S',\tau^*_S\pi_*\cG)\isom\Gamma(S',\pi'_*\tau^*_V\cG)
\isom\Gamma(V',\tau^*_V\cG).
$$
Taking $\cG:=\psi_*\cO_{\!X_\et}$, we then obtain
the sought isomorphism $\omega_\emptyset$
$$
A=\Gamma(X,\cO_{\!X_\et})\isom\Gamma(V,\psi_*\cO_{\!X_\et})
\isom\Gamma(V',\cF).
$$
Likewise, the map $\omega_{0,\emptyset}$ is defined as the composition
$$
A_0=\Gamma(S,\cO_{\!S})\to\Gamma(S,\pi_*\cO_V)\isom
\Gamma(V',\tau^*_V\cO_V)
$$
from which the commutativity of $\cD_\emptyset$ is straightforward.
Next, let $\Lambda\neq\emptyset$, so that $V_\Lambda$ is an
affine scheme; set $V^\he_\Lambda:=\Spec\,A_{0,\Lambda}^\he$ and
$X_\Lambda:=X\cap\Spec\,A_\Lambda$. We let
$\psi_\Lambda:X_\Lambda\to V_\Lambda$ be the restriction of
$\psi$, and we consider the cartesian diagram of schemes
$$
\xymatrix{ &
X^\he_\Lambda \ar[r]^-{\tau''_{X_\Lambda}} \ar[d]_{\psi^\he_\Lambda} &
X_\Lambda \ar[d]^{\psi_\Lambda} \\
V'_\Lambda \ar[r]^-{\tau'_{V_\Lambda}} &
V^\he_\Lambda \ar[r]^-{\tau''_{V_\Lambda}} & V_\Lambda
}$$
where $\tau'_{V_\Lambda}$ and $\tau''_{V_\Lambda}$ are induced by
the natural projection
$A^\he_{0,\Lambda}\to A^\he_{0,\Lambda}/I_{0,\Lambda}^\he\isom
A_{0,\Lambda}/I_0A_{0,\Lambda}$, and respectively the henselization
map $A_{0,\Lambda}\to A^\he_{0,\Lambda}$. Then $\omega_\Lambda$
is defined as the composition of isomorphisms of $A$-algebras
$$
\begin{aligned}
A^\he_\Lambda=\Gamma(X^\he_\Lambda,\cO_{\!X^\he_{\Lambda,\et}})\isom
\Gamma(V^\he_\Lambda,
\psi^\he_{\Lambda*}\circ\tau''^*_{X_\Lambda}\cO_{\!X_{\Lambda,\et}})\isom\,&
\Gamma(V^\he_\Lambda,
\tau''^*_{V_\Lambda}\circ\psi_{\Lambda*}\cO_{\!X_{\Lambda,\et}}) \\
\isom\,&
\Gamma(V'_\Lambda,\tau'^*_{V_\Lambda}\circ\tau''^*_{V_\Lambda}\circ
\psi_{\Lambda*}\cO_{\!X_{\Lambda,\et}}) \\
\isom\,& \Gamma(V'_\Lambda,\cF)
\end{aligned}
$$
where the bijectivity of the third map is due to
\cite[Exp.XII, Prop.6.5(i)]{SGA4-3}, and that of the 
second one is due to the fact that $\tau''_{V_\Lambda}$ is
the limit of a cofiltered system of affine \'etale
$V_\Lambda$-schemes. Likewise, $\omega_{0,\Lambda}$ is
obtained as the composition
$$
A^\he_{0,\Lambda}=\Gamma(V^\he_\Lambda,\cO_{V^\he_\Lambda})\isom
\Gamma(V'_\Lambda,\tau'^*_{V_\Lambda}\cO_{V^\he_\Lambda})\isom
\Gamma(V'_\Lambda,\tau^*_V\cO_V).
$$
Then the commutativity of $\cD_\Lambda$ follows by a
direct inspection.

(ii): We consider the diagram of $A_0$-algebras
$$
\xymatrix{
A^\he_{0,\Lambda'} \ar[r]^-{\omega_{0,\Lambda'}}
\ar[d]_{\mu^\he_{0,\Lambda',\Lambda}} &
\Gamma(V'_{\Lambda'},\tau_V^*\cO_{V,\et}) \ar[r] \ar[d] &
\Gamma(V'_{\Lambda'},\cF) \ar[d]_{\rho_{\Lambda',\Lambda}} &
\Gamma(V_{\Lambda'},\psi_*\cO_{\!X,\et})=A_{\Lambda'}
\ar[l] \ar[d]^{\mu_{\Lambda',\Lambda}} \\
A^\he_{0,\Lambda} \ar[r]^-{\omega_{0,\Lambda'}} &
\Gamma(V'_\Lambda,\tau_V^*\cO_{V,\et}) \ar[r] &
\Gamma(V'_\Lambda,\cF) &
\Gamma(V_\Lambda,\psi_*\cO_{\!X,\et})=A_\Lambda \ar[l]
}$$
whose central square subdiagram is induced by the natural
morphism $\cO_{V,\et}\to\psi_*\cO_{\!X,\et}$, and the right
square subdiagram is induced by the unit of adjunction
$\psi_*\cO_{\!X,\et}\to\tau_{V*}\tau^*_V\psi_*\cO_{\!X,\et}$
for the adjoint pair $(\tau^*_V,\tau_{V*})$, and where
$\mu_{\Lambda',\Lambda}$ and $\mu^\he_{0,\Lambda',\Lambda}$ are
defined as in remark \ref{rem_extract-proj}(ii). Therefore,
both these subdiagrams commute; the left square subdiagram
commutes as well, due to the universal property of
henselization; taking into account (i), the assertion
follows (details left to the reader).
\end{pfclaim}

Lastly, we have a natural map
$$
A\xrightarrow{\ \alpha\ }\Equal\Bigl(\prod_{i=0}^n
\xymatrix{A^\he_{\{i\}} \ar@<-.5ex>[r]_-\beta \ar@<.5ex>[r]^-\gamma &}
\prod_{i,j=0}^n A^\he_{\{i,j\}}\Bigr)
$$
induced by the maps $\mu^\he_{\emptyset,\{i\}}$, where $\beta$
and $\gamma$ are deduced from the maps $\mu^\he_{\{i\},\{i,j\}}$.
From claim \ref{cl_I-am-convinced} we deduce that $\alpha$
is an isomorphism, and taking into account remark
\ref{rem_extract-proj}(iii), the theorem follows.
\end{proof}

\begin{remark}\label{rem_extend-to-whole-top}
Let $\underline X$ be any topologically local quasi-affinoid
scheme.

(i)\ \
As explained in \cite[Ch.0, \S3.2.2]{EGAI}, every sheaf on
$\cQ(\underline X)$ admits a natural extension to a sheaf
$\cF'$ on the topological space $\Spa\,\underline X$. Namely,
if $U\subset\Spa\,\underline X$ is any open subset, one sets
$$
\cF'(U):=\lim_{V\in(\cQ(\underline X)/U)^o}\cF(V)
$$
where $\cQ(\underline X)/U$ denotes the full subcategory
of $\cQ(\underline X)$ whose objects are the quasi-affinoid
open subsets of $\Spa\,\underline X$ contained in $U$. We
shall use the notation $\cO^\he_{\Spa\,\underline X}$ and
$\cO^{\he+}_{\Spa\,\underline X}$ also to refer to the sheaves
on $\Spa\,\underline X$ that extend naturally the respective
sheaves on $\cQ(\underline X)$.

(i)\ \
Simple examples show that, in most cases, the presheaves
$\cO^\loc_{\Spa\,\underline X}$ and $\cO^{\loc+}_{\Spa\,\underline X}$
are {\em not} sheaves on the site $\cQ(\underline X)$. On
the other hand, , in several situations of interest the
presheaves $\cO^\wedge_{\Spa\,\underline X}$ and
$\cO^{\wedge+}_{\Spa\,\underline X}$ are sheaves of topological
rings. In Huber's work \cite{Hu1}, this sheaf property is
proven under two different types of assumptions,
corresponding roughly to the two main classes of
non-archimedean analytic spaces that are encountered
in applications : the ``generic fibers'' of locally
noetherian formal schemes, and the rigid analytic
varieties locally of finite type over a complete rank
one valued field. In the literature, one finds also
some attempts to unify these cases (and the few other
known ones) under a single axiomatic framework; hereafter
we shall adopt the approach (though, not the terminology)
proposed in \cite{Fu-Ga-Ka}, that is based on the notions
introduced in section \ref{sec_anal-noetherian}.
Later, we shall prove the sheaf property also in the case
where $\cO_{\!X}(X)$ admits a perfectoid ring of definition :
this is a completely different situation, for which we shall
need to develop an {\em ad hoc} method.
\end{remark}

\begin{definition}
(i)\ \ 
Let $A$ be any f-adic ring. We say that $A$ is
{\em f-adic analytically noetherian} (resp. {\em f-adic
universally analytically noetherian}), if it admits an
analytically noetherian (resp. universally analytically
noetherian) ring of definition (see definition
\ref{def_anal-noetherian}(i)).

(ii)\ \ 
Let $\underline A:=(A,A^+)$ be any affinoid ring. We say that
$\underline A$ is {\em analytically noetherian} (resp. {\em
universally analytically noetherian}) if $A$ is f-adic
analytically noetherian (resp. f-adic universally analytically
noetherian).

(iii)\ \
Let $\underline A:=(A,A^+,U)$ be any quasi-affinoid ring.
We say that $\underline A$ is {\em analytically noetherian}
(resp. {\em universally analytically noetherian}) if the
same holds for $(A,A^+)$.

(iv)\ \
Let $\underline X$ be any quasi-affinoid scheme. We say
that $\underline X$ is {\em analytically noetherian}
(resp. {\em universally analytically noetherian}) if the
same holds for $\sGamma(\underline X)$.
\end{definition}

\sset\subsubsection{}\label{subsec_f-adic-univ-an-noether}
Let $A$ be any f-adic ring, and $n\in\N$ any integer.
According to proposition \ref{prop_polynomial-top-rings}(i,ii),
there exists a unique f-adic topology $\cT_n$ on $A[T_1,\dots,T_n]$
such that, if $A_0$ is any subring of definition of $A$,
then $A_0[T_1,\dots,T_n]$ is a subring of definition of
$A[T_1,\dots,T_n]$, and if $I\subset A_0$ is an ideal of
adic definition, then $I[T_1,\dots,T_n]$ is an ideal of
adic definition for $A_0[T_1,\dots,T_n]$.

\begin{lemma}\label{lem_f-adic-noether}
With the notation of \eqref{subsec_f-adic-univ-an-noether},
the following holds :
\begin{enumerate}
\item
$A$ is f-adic analytically noetherian (resp. f-adic universally
analytically noetherian) if and only if every ring of definition
of $A$ is analytically noetherian (resp. universally analytically
noetherian).
\item
$A$ is f-adic universally analytically noetherian if and only
if $(A[T_1,\dots,T_n],\cT_n)$ is f-adic analytically noetherian
for every $n\in\N$.
\item
Let $B\subset A$ be any open subring, and endow $B$ with the
topology induced from $A$. Then $A$ is f-adic analytically
noetherian (resp. f-adic universally analytically noetherian)
if and only if the same holds for $B$.
\item
Let $f:A\to A'$ be a surjective ring homomorphism, and endow
$A'$ with the topology induced by $A$ via $f$. Then, if $A$
is f-adic analytically noetherian (resp. f-adic universally
analytically noetherian), the same holds for $A'$.
\item
If $\underline A$ is any f-adic analytically noetherian (resp.
f-adic universally analytically noetherian) quasi-affinoid
ring, $\Spec\,\underline A$ is an analytically noetherian
(resp. universally analytically noetherian) quasi-affinoid
scheme.
\item
Let $\underline X$ be a topologically local and universally
analytically noetherian quasi-affinoid scheme, and
$\underline Y\to\underline X$ a morphism of topologically
local quasi-affinoid schemes that represents a rational
subset of\/ $\Spa\,\underline X$. Then $\underline Y$ is
universally analytically noetherian.
\end{enumerate}
\end{lemma}
\begin{proof}(i): Let $A_0,A_1\subset A$ be two rings of
definition, and let us show that $A_0$ is analytically
noetherian if and only if the same holds for $A_1$. To
this aim, suppose first that $A_0\subset A_1$, and let
$I\subset A_0$ be any finitely generated ideal of adic
definition, so that $IA_1$ is an ideal of adic definition
for $A_1$. Lemma \ref{lem_deja-vu}(iii) says that the
inclusion map $A_0\to A_1$ identifies the analytic locus
$U_0$ of $\Spec\,A_0$ with the analytic locus $U_1$ of
$\Spec\,A_1$, so $U_0$ is a noetherian scheme if and only
if the same holds for $U_1$. Next, suppose that $A_0$ is
analytically noetherian; let $M_1$ be any $A_1$-module of
finite type, $x_1,\dots,x_n$ any system of generators for
$M_1$, and set $M_0:=A_0x_1+\cdots+A_0x_n\subset M_1$. By
proposition \ref{prop_f-adics}(ii), there exists $n\in\N$
such that $I^nA_1\subset A_0$, whence $I^nM_1\subset M_0$.
On the other hand, by assumption there exists $r\in\N$ such
that $\Ann_{M_0}(I^k)=\Ann_{M_0}(I^r)$ for every $k\geq r$.
Thus, let $x\in\Ann_{M_1}(I^kA_1)$ for some $k\geq r$; we
get $I^nx\in\Ann_{M_1}(I^kA_1)\cap M_0=\Ann_{M_0}(I^k)$, so
$I^{n+r}x=0$, as asserted, in this case.

Suppose next that $A_1$ is analytically noetherian,
and let $M_0$ be any $A_0$-module of finite type; the
short exact sequence of $A_0$-modules
$0\to A_0\to A_1\to A_1/A_0\to 0$ induces an exact sequence
$$
T_1:=\Tor_1^{A_0}(A_1/A_0,M_0)\to M_0
\xrightarrow{\ \phi\ }M_1:=A_1\otimes_{A_0}M_0
$$
and notice that $I^nT_1=0$. By assumption, there exists
$r\in\N$ such that $\Ann_{M_1}(I^kA_1)=\Ann_{M_1}(I^rA_1)$
for every $k\geq r$. Suppose then that $x\in\Ann_{M_0}(I^k)$
for some $k\geq r$; it follows easily that
$I^rx\subset\Ker\,\phi$, whence $I^{r+n}x=0$, as required.

Lastly, if $A_0$ and $A_1$ are arbitrary rings of definition
for $A$, the same holds for $A_2:=A_0\cdot A_1$ (corollary
\ref{cor_f-adics}(i)). Suppose that $A_0$ is analytically
noetherian; by the foregoing, it follows that the same holds
for $A_2$, and then by the same token it holds for $A_1$
as well.

The assertion for the universally noetherian case is an
immediate consequence.

Assertion (ii) follows directly from (i) : details
left to the reader.

(iii) follows from (i) and corollary \ref{cor_f-adics}(ii).

(iv): First, we know that $A'$ is f-adic, by example
\ref{ex_f-adic-quotient}(iii), and if $A_0\subset A$ is
a ring of definition, $f(A_0)$ is a ring of definition
of $A'$. But if $A_0$ is analytically noetherian (resp.
universally analytically noetherian), the same holds for
$f(A_0)$, by lemma \ref{lem_sorite-analyt}(ii). The
assertion follows.

(v): Say that $\underline A=(A,A^+,U)$, and let
$\rho:A\to A_U:=\cO_{\Spec\,A}(U)$ be the restriction map.
By definition, the topology $\cT_{\rho(A)}$ on $\rho(A)$
induced by the inclusion into $A_U$ agrees with the quotient
topology induced by the surjection $A\to\rho(A)$. From (iv)
we deduce that $(\rho(A),\cT_{\rho(A)})$ is f-adic analytically
noetherian (resp. f-adic universally analytically noetherian).
Then, the assertion follows from (iii).

(vi) follows easily from (v) and lemma \ref{lem_sorite-analyt}(ii),
after inspecting the proof of theorem \ref{th_represent-rational}.
\end{proof}

\begin{lemma}\label{lem_simple-completion}
Let $\underline X:=(X,\cT_X,A^+_X)$ be any analytically
noetherian quasi-affinoid scheme. Then the natural morphism
$\sGamma(\underline X)^\wedge\to\sGamma(\underline X^\wedge)$
is an isomorphism (see remark
{\em\ref{rem_depth-and-completion}(ii)}).
\end{lemma}
\begin{proof} Let $A^\wedge_X$ be the separated completion of
$A_X:=\cO_{\!X}(X)$, and set $Y:=\Spec\,A^\wedge_X\times_{\Spec\,A_X}X$;
in light of proposition \ref{prop_like-noetherian}(ii) and
corollary \ref{cor_base-change-where}, the natural map
$A^\wedge_X\to\cO_Y(Y)$ is an isomorphism. Then the assertion
follows from remark \ref{rem_depth-and-completion}(ii).
\end{proof}

\sset\subsubsection{}\label{subsec_smaller-rational-site}
For any topologically local quasi-affinoid scheme $\underline X$,
we may consider the site
$$
(\cR(\underline X),J_\cR)
$$
where $\cR(\underline X)$ is the full subcategory of
$\cQ(\underline X)$ whose objects are the rational subsets
of $\Spa\,\underline X$. The sieves covering a given object
$R$ of $\cR(\underline X)$ for the topology $J_\cR$ are those
generated by the families
$R_\bullet:=(R_\lambda~|~\lambda\in\Lambda)$ of rational subsets of
$\Spa\,\underline X$ such that $\bigcup_{\lambda\in\Lambda}R_\lambda=R$.
For every such $R_\bullet$, after fixing a total ordering on
$\Lambda$ we get an augmented alternating \v{C}ech complex
$C^\bullet_\mathrm{alt}(R_\bullet,\cO^\wedge_{\Spa\,\underline X})$, whose
degree $n$ term, for every $n\geq 0$, is a product of topological
rings of the form
$\cO^\wedge_{\Spa\,\underline X}(R_{\lambda_0}\cap\cdots\cap R_{\lambda_n})$,
with $\lambda_0<\cdots<\lambda_n$ ranging over all strictly
increasing sequences of elements of $\Lambda$ (see remark
\ref{rem_justify-name}(i).
We endow $C^n_\mathrm{alt}(R_\bullet,\cO^\wedge_{\Spa\,\underline X})$
with the corresponding product topology, and we notice that
this topology is independent of the chosen ordering for
$\Lambda$, and the differentials $d^\bullet$ of the \v{C}ech
complex are continuous maps. With this notation, we have:

\begin{theorem}\label{th_compl-an-noeth-are-sheaves}
For any topologically local and universally analytically
noetherian quasi-affinoid scheme
$\underline X:=(X,\cT_X,A^+_X)$, and every
$R\in\Ob(\cR(\underline X))$, the following holds :
\begin{enumerate}
\item
The presheaves $\cO^{\wedge+}_{\Spa\,\underline X}$ and
$\cO^\wedge_{\Spa\,\underline X}$ are sheaves of topological rings
on $(\cR(\underline X),J_\cR)$.
\item
For every finite covering $\fU$ of $R$ consisting of rational
subsets, the differentials of
$C^\bullet_\mathrm{alt}(\fU,\cO^\wedge_{\Spa\,\underline X})$
are strict and its cohomology has the discrete topology.
\end{enumerate}
\end{theorem}
\begin{proof}(i): Arguing as in the proof of theorem
\ref{th_henselians-are-sheaves}, we reduce to checking the
assertion for $\cO^\wedge_{\Spa\,\underline X}$. We remark :
\begin{claim}\label{cl_tidy-up}
For any given $R\in\Ob(\cR(\underline X))$, let
$\phi_{\underline Y/\underline X}:\underline Y\to\underline X$
be any morphism of quasi-affinoid schemes representing the
sub-presheaf $h_R$ of $h_{\underline X}$; then the morphism
$\phi_{\underline Y/\underline X}$ induces an equivalence of
categories
$$
\cR(\underline Y)\isom\cR(\underline X)/R
\qquad
V\mapsto\Spa\,\phi_{\underline Y/\underline X}(V)\subset R.
$$
Moreover, $\underline Y$ is still universally analytically
noetherian.
\end{claim}
\begin{pfclaim} These assertions follow from lemmata
\ref{lem_rat-in-rat} and \ref{lem_f-adic-noether}(vi).
\end{pfclaim}

Now, say that $\sGamma(\underline X)=(A,A^+,X)$ and let
$R_\bullet:=(R_0,\dots,R_n)$ be the standard covering of
$\Spa\,\underline X$ attached to a given sequence
$f_\bullet:=(f_0,\dots,f_n)$ of elements of $A$; define also
the rational subset $R_\Lambda$ as in remark
\ref{rem_extract-proj}(i), for every subset
$\Lambda\subset\{0,\dots,n\}$. In light of claim
\ref{cl_tidy-up} and remark \ref{rem_sheaves-with-values-in-A}(i),
we are reduced to showing that the natural map
\set\begin{equation}\label{eq_global-to-equal}
\cO^\wedge_{\Spa\,\underline X}(\Spa\,\underline X)\to
\Equal\Bigl(\prod_{i=0}^n\xymatrix{\cO^\wedge_{\Spa\,\underline X}(R_i)
\ar@<-.5ex>[r] \ar@<.5ex>[r] &}
\prod_{i,j=0}^n\cO^\wedge_{\Spa\,\underline X}(R_{\{i,j\}})\Bigr)
\end{equation}
is an isomorphism of topological rings. To this aim, we
pick a subring of definition $A_0\subset A$, a finitely
generated ideal $I_0\subset A_0$ of adic definition, we
set $S:=\Spec\,A_0$, we denote by $J\subset A$ the ideal
generated by $f_\bullet$, and we consider the system
$(\underline A{}_\Lambda~|~\Lambda\subset\{0,\dots,n\})$
of quasi-affinoid rings, and their rings of definition
$A_{0,\Lambda}\subset A_\Lambda$ constructed in remark
\ref{rem_extract-proj}(i). Let also $A^\wedge_\Lambda$ be
the separated completion of $A_\Lambda$, and notice that the map
$$
\omega^\wedge_\Lambda:A_\Lambda^\wedge\to
\cO^\wedge_{\Spa\,\underline X}(R_\Lambda)
$$
of remark \ref{rem_extract-proj}(ii) is an isomorphism of
topological $A$-algebras for every $\Lambda\subset\{0,\dots,n\}$,
by lemma \ref{lem_simple-completion}. Define $S$, $V$,
$\psi:X\to V$ and $V_\Lambda$ for every subset
$\Lambda\subset\{0,\dots,n\}$ as in remark
\ref{rem_extract-proj}(v); set also $S':=\Spec\,A_0/I_0$
and let $V^\wedge$ be the completion of $V$ along $S'\times_SV$,
and $\pi:V^\wedge\to V$ the induced morphism of locally ringed
spaces. We define
$$
\cF:=\pi^*\psi_*\cO_{\!X}
\qquad\text{and}\qquad
V^\wedge_\Lambda:=\pi^{-1}V_\Lambda
\qquad\text{for every $\Lambda\subset\{0,\dots,n\}$}
$$
and we regard $\cF$ as a presheaf of topological abelian
groups, as follows. The natural morphism $\cO_V\to\psi_*\cO_{\!X}$
induces a morphism
$$
\phi:\cO_{V^\wedge}\to\cF
$$
and for every open subset $U\subset V^\wedge$ we endow
$\cF(U)$ with the unique group topology such that the
resulting map $\phi_U:\cO_{V^\wedge}(U)\to\cF(U)$ is continuous
and open. It follows easily that for every inclusion
$U'\subset U$ of open subsets of $V^\wedge$, the restriction
map $\cF(U)\to\cF(U')$ is continuous.

\begin{claim}\label{cl_I-am-completely-convinced}
(i)\ \ 
$\cF(V^\wedge_\Lambda)$ is a topological ring for every
$\Lambda\subset\{0,\dots,n\}$, and there exists a natural
isomorphism of topological $A$-algebras
$$
\omega_\Lambda:A^\wedge_\Lambda\isom\cF(V^\wedge_\Lambda).
$$
\begin{enumerate}
\addenu
\item
For every $\Lambda'\subset\Lambda$, let
$\rho_{\Lambda',\Lambda}:\cF(V^\wedge_{\Lambda'})\to
\cF(V^\wedge_\Lambda)$ be the restriction homomorphism
associated with the inclusion
$V^\wedge_\Lambda\subset V^\wedge_{\Lambda'}$. Then we have a
commutative diagram
$$
\xymatrix{ A^\wedge_{\Lambda'} \ar[r]^-{\omega_{\Lambda'}}
\ar[d]_{\mu^\wedge_{\Lambda',\Lambda}} & \cF(V^\wedge_{\Lambda'})
\ar[d]^{\rho_{\Lambda',\Lambda}} \\
A^\wedge_\Lambda \ar[r]^-{\omega_\Lambda} &
\cF(V^\wedge_\Lambda).
}$$
\end{enumerate}
\end{claim}
\begin{pfclaim}(i): We consider first the case where
$\Lambda=\emptyset$, and recall that $V_\emptyset=V$.
Notice also that $\psi_*\cO_{\!X}$ is a quasi-coherent
$\cO_V$-module, by virtue of \cite[Ch.I, Cor.9.2.2]{EGAI}.
Since $V$ is a projective $S$-scheme, corollary
\ref{cor_rompiballe}(iii) and remark \ref{rem_conclusion} yield
a natural isomorphism
$$
\omega_\emptyset:A^\wedge\isom A^\wedge_0\otimes_{A_0}A
\isom A^\wedge_0\otimes_{A_0}\Gamma(V,\psi_*\cO_{\!X})
\isom\cF(V^\wedge)
$$
fitting into a commutative diagram
$$
\xymatrix{
A^\wedge_0 \ar[rr] \ar[d] & & A^\wedge \ar[d]^{\omega_\emptyset} \\
\cO_{V^\wedge}(V^\wedge) \ar[rr]^-{\phi_{V^\wedge}} & &
\cF(V^\wedge)
}$$
whose left vertical arrow is the $I_0$-adic completion of
the natural map $A_0\to\cO_V(V)$, and whose top horizontal
arrow is the inclusion map. It follows already that
$\omega_\emptyset$ is a continuous map, and it thus remains
only to check that $\omega_\emptyset$ is an open map.
To this aim, it suffices to check that the same holds for
the left vertical arrow $A^\wedge_0\to\cO_{V^\wedge}(V^\wedge)$.
However, $\cO_V(V)$ is an $A_0$-module of analytically
finite type, by theorem \ref{th_analyt-proper-finiteness},
and the support of $\cO_V(V)/A_0$ lies in the non-analytic
locus of $S$, by virtue of remark \ref{rem_extract-proj}(vi);
by  remark \ref{rem_anal-noetherian}(iii) we deduce that
$I^t_0\cdot\cO_V(V)\subset A_0$ for some $t\in\N$, whence
$I^t_0\cdot\cO_{V^\wedge}(V^\wedge)\subset A_0^\wedge$. Since the
topology of $\cO_{V^\wedge}(V^\wedge)$ is $I_0$-adic (corollary
\ref{cor_rompiballe}(ii)), the assertion follows.

In case $\Lambda$ is not empty, $V_\Lambda$ is affine, and
a direct inspection yields a natural isomorphism of topological
rings $A^\wedge_{0,\Lambda}\isom\cO_{V^\wedge}(V^\wedge_\Lambda)$.
We may then appeal to proposition \ref{prop_up-to-completion}(iii)
to see that the natural map
\set\begin{equation}\label{eq_affine-case}
A^\wedge_\Lambda\isom A^\wedge_{0,\Lambda}\otimes_{A_{0,\Lambda}}A_\Lambda
\isom A^\wedge_{0,\Lambda}\otimes_{A_{0,\Lambda}}
\Gamma(V_\Lambda,\psi_*\cO_{\!X})\to\cF(V^\wedge_\Lambda)
\end{equation}
is also an isomorphism of rings, fitting into another
commutative diagram
\set\begin{equation}\label{eq_noise}
{\diagram
A^\wedge_{0,\Lambda} \ar[rr] \ar[d] & & A^\wedge_\Lambda  \ar[d] \\
\cO_{V^\wedge}(V^\wedge_\Lambda) \ar[rr]^-{\phi_{V^\wedge_\Lambda}} & &
\cF(V^\wedge_\Lambda)
\enddiagram}
\end{equation}
whose top horizontal arrow is again the inclusion map,
so $\cF(V^\wedge_\Lambda)$ is also a topological ring with the
topology defined in the foregoing, and \eqref{eq_affine-case}
is even an isomorphism of topological rings, if
$\Lambda\neq\emptyset$.

(ii): The abelian groups appearing in the diagram are all
endowed with natural $A^\wedge_{\Lambda'}$-module structures,
and it is easily seen that all the arrows are both
$A^\wedge_{0,\Lambda'}$-linear and $A_{\Lambda'}$-linear maps,
whence the assertion.
\end{pfclaim}

Next, consider the affine open covering
$V^\wedge_\bullet:=(V^\wedge_{\{i\}}~|~i=0,\dots,n)$ of $V^\wedge$,
and endow the terms of the augmented alternating \v{C}ech
complex $C^\bullet_\mathrm{alt}(V^\wedge_\bullet,\cF)$ with the
product topologies, as in \eqref{subsec_smaller-rational-site}.
Combining claim \ref{cl_I-am-completely-convinced} and
remark \ref{rem_extract-proj}(iii) we deduce an isomorphism
of complexes of topological abelian groups
\set\begin{equation}\label{eq_identify-Cech}
C^\bullet_\mathrm{alt}(V^\wedge_\bullet,\cF)\isom
C^\bullet_\mathrm{alt}(R_\bullet,\cO^\wedge_{\Spa\,\underline X}).
\end{equation}
Hence, to conclude it suffices to show :

\begin{claim} The complex
$C^\bullet_\mathrm{alt}(V^\wedge_\bullet,\cF)$ has strict
differentials, and its cohomology has the discrete topology.
\end{claim}
\begin{pfclaim} Diagram \eqref{eq_noise} implies that
$\phi_{V^\wedge_\Lambda}$ is injective for every $\Lambda\neq\emptyset$;
it follows easily that the same holds for $\phi_{V^\wedge}$.
Hence the induced map of alternating \v{C}ech complexes
$$
C^\bullet_\mathrm{alt}(V^\wedge_\bullet,\phi):
C^\bullet_\mathrm{alt}(V^\wedge_\bullet,\cO_{V^\wedge})\to
C^\bullet_\mathrm{alt}(V^\wedge_\bullet,\cF)
$$
is injective in every degree, and by construction
$C^i_\mathrm{alt}(V^\wedge_\bullet,\phi)$ is an open map for
every $i\in\Z$, so we are reduced to showing that the
differentials $d^i$ of
$C^\bullet_\mathrm{alt}(V^\wedge_\bullet,\cO_{V^\wedge})$ are
strict and induce open maps
$C^i_\mathrm{alt}(V^\wedge_\bullet,\cO_{V^\wedge})\to\Ker\,d^{i+1}$
for every $i\in\Z$. If $i<-1$,
there is nothing to show. For $i=-1$, the assertion follows
from remark \ref{rem_sheaves-with-values-in-A}(i). In the remaining
cases, the assertion comes down to the following. For every
$i,t\in\N$ there exists $s\in\N$ such that
\set\begin{equation}\label{eq_kis}
I^s_0\cdot C^{i+1}_\mathrm{alt}(V^\wedge_\bullet,\cO_{V^\wedge})\cap
\Ker\,d^{i+1}\subset I^t_0\cdot\Img\,d^i.
\end{equation}
For every $k\in\N$, denote by $d^\bullet_k$ the differentials
of the \v{C}ech complex
$C^\bullet_\mathrm{alt}(V^\wedge_\bullet,I^k_0\cO_{V^\wedge})$;
condition \eqref{eq_kis} is equivalent to
$$
\Ker\,d^{i+1}_s\subset\Img\,d^i_t.
$$
Taking into account theorem \ref{th_coh-vanish}(i), we are
therefore reduced to showing that for every $i,t\in\N$
there exists an integer $s\geq t$ such that the natural map
$$
H^{i+1}(V^\wedge,I^s_0\cO_{V^\wedge})\to H^{i+1}(V^\wedge,I^t_0\cO_{V^\wedge})
$$
vanishes. Taking into account corollaries
\ref{cor_analyt-finiteness}(i) and \ref{cor_rompiballe}(i),
we are further reduced to showing that for every $i,t\in\N$
there exists $s\in\N$ such that $I^s\cdot H^{i+1}(V,I^t_0\cO_V)=0$.
Notice that the $A_0$-module $H^{i+1}(V,I^t_0\cO_V)$ is analytically
of finite type (theorem \ref{th_analyt-proper-finiteness} and
remark \ref{rem_conclusion}(i)); in light of remark
\ref{rem_anal-noetherian}(iii) it then suffices to check that
the support of $H^{i+1}(V,I^t_0\cO_V)$ lies in the non-analytic
locus of $S$. The latter follows easily from remark
\ref{rem_extract-proj}(vi) : details left to the reader.
\end{pfclaim}

(ii): Arguing as in the proof of (i), we reduce easily to the
case where $R=\Spa\,\underline X$. Say that $\fU=(U_i~|~i\in I)$
for some set $I$, and pick any standard covering
$\fU':=(U'_i~|~i\in I')$ of $\Spa\,\underline X$ that refines
$\fU$; after fixing total orderings for $I$ and $I'$ we may
consider the double complex
$$
C_\alt^{\bullet\bullet}(\fU,\fU'):=
\Hom^\bullet_{\Z_X}(R^\alt_{\bullet\bullet}(\fU,\fU'),\cO^\wedge_{\Spa\,\underline X})
$$
(notation of remark \ref{rem_double-cechs}(ii)) which vanishes
in every bidegree $(p,q)$ with either $p<-1$ or $q<-1$. For
every $n\in\N$, define $I^{n+1}_\alt$ and $I'^{n+1}_\alt$ as in
\eqref{subsec_choose-total}, and for every
$\underline t\in I^{n+1}_\alt$ and every
$\underline t'\in I'^{n+1}_\alt$, set
$\fU'_{\underline t}:=(U'_{i'}\cap U_{\underline t}~|~i'\in I')$ and
$\fU_{\underline t'}:=(U_i\cap U_{\underline t'}~|~i\in I)$. We notice
that :
\begin{itemize}
\item
The complex $C_\alt^{\bullet,-1}(\fU,\fU')$ coincides with
$C_\alt^\bullet(\fU,\cO^\wedge_{\Spa\,\underline X})$.
\item
The complex $C_\alt^{-1,\bullet}(\fU,\fU')$ coincides with
$C_\alt^\bullet(\fU',\cO^\wedge_{\Spa\,\underline X})$.
\item
For every $n\in\N$, the complex $C_\alt^{\bullet,n}(\fU,\fU')$
is the product $\prod_{\underline t'\in I'^{n+1}_\alt}
C^\bullet_\alt(\fU_{\underline t'},\cO^\wedge_{\Spa\,\underline X})$.
\item
For every $n\in\N$, the complex $C_\alt^{n,\bullet}(\fU,\fU')$
is the product $\prod_{\underline t\in I^{n+1}_\alt}
C^\bullet_\alt(\fU'_{\underline t},\cO^\wedge_{\Spa\,\underline X})$.
\end{itemize}
Since $\fU'_{\underline t}$ is a standard covering of $U_{\underline t}$,
the proof of (i) shows that the differentials of the complex
$C_\alt^\bullet(\fU'_{\underline t},\cO^\wedge_{\Spa\,\underline X})$ induce
open maps $\delta^{n,\bullet}_v$ as in \eqref{eq_strict-cont-double},
for every $n\in\N$. Likewise, the same holds for the induced maps
$\delta^{-1,\bullet}_v$. Moreover, since $\fU'$ refines $\fU$, for
every $n\in\N$ and every $\underline t'\in I'^{n+1}_\alt$ there
exists $i\in I$ such that $U'_{\underline t'}\subset U_i$, in which
case remark \ref{rem_explicit-homotopy} yields a homotopy $h^\bullet$
from the identity automorphism of
$C^\bullet_\alt(\fU_{\underline t'},\cO^\wedge_{\Spa\,\underline X})$ to the
zero map, and a simple inspection shows that $h^i$ is a continuous
map, for every $i\in\Z$. Then, lemma \ref{lem_contract-cont}(i)
says that the differentials of
$C^\bullet_\alt(\fU_{\underline t'},\cO^\wedge_{\Spa\,\underline X})$ induce
open maps $\delta^{\bullet,n}_h$ as well, for every $n\in\N$. Lastly,
we invoke corollary \ref{cor_double-strict-continuous} to conclude.
\end{proof}

\sset\subsubsection{}\label{subsec_support-for-qaff}
Let $\underline X:=(X,\cT_{\!X},A^+_X)$ be a quasi-affinoid
scheme fulfilling the assumptions of theorem
\ref{th_compl-an-noeth-are-sheaves}, so that
$\cO^\wedge_{\Spa\,\underline X}$ is a sheaf of topological rings
on the site $(\cR(\underline X),J_\cR)$.
Following \cite[Ch.0, \S3.2.2]{EGAI} (and see also remark
\ref{rem_extend-to-whole-top}(i)), we may extend (uniquely
up to unique isomorphism) both $\cO^\wedge_{\Spa\,\underline X}$
and $\cO^{\wedge+}_{\Spa\,\underline X}$ to sheaves of topological
groups on the topological space $\Spa\,\underline X$, and we
shall denote these extensions with the same names. Then the
pair $(\Spa\,\underline X,\cO^\wedge_{\Spa\,\underline X})$ is a
topologically ringed space. Moreover, lemma
\ref{lem_stalks-are-local} says that $\Spa\,\underline X$
is a locally ringed space, so there exists a unique morphism
of locally ringed spaces
\set\begin{equation}\label{eq_support-map-for-X}
\sigma_{\!\underline X}:(\Spa\,\underline X,\cO^\wedge_{\Spa\,\underline X})
\to\Spec\,A
\end{equation}
such that the corresponding map
$\cO_{\Spec\,A}\to\sigma_{\underline X*}\cO^\wedge_{\Spa\,\underline X}$
yields on global sections the natural map
$A\to A^\wedge:=\cO^\wedge_{\Spa\,\underline X}(\Spa\,\underline X)$
(\cite[Ch.I, Prop.1.6.3]{EGAI-new}). By inspecting {\em loc.cit.}
it is easily seen that $\sigma_{\!\underline X}$ is the
composition
$$
\Spa\,\underline X\subset\Cont\,A
\xrightarrow{\ \sigma_{\!A}\ }\Spec\,A
$$
where $\sigma_{\!A}$ is the restriction of the support map
(see remark \ref{rem_Spv-of-ring}(iii)). Especially, the
image of $\sigma_{\!\underline X}$ lies in $X$, whence a well
defined morphism of locally ringed spaces
$$
(\Spa\,\underline X,\cO^\wedge_{\Spa\,\underline X})
\to(X,\cO_{\!X}).
$$

\begin{corollary}\label{cor_support-for-qaff}
With the notation of \eqref{subsec_support-for-qaff}, the
morphism $\sigma_{\underline X}$ induces an isomorphism
$$
A^\wedge\otimes_AH^i(X,\cO_{\!X})\isom
H^i(\Spa\,\underline X,\cO^\wedge_{\Spa\,\underline X})
\qquad
\text{for every $i\in\N$}.
$$
\end{corollary}
\begin{proof} Let us show first the following special case :
\begin{claim}\label{cl_loud}
In the situation of the corollary, suppose moreover that
$\underline X$ is affinoid. Then
$H^i(\Spa\,\underline X,\cO^\wedge_{\Spa\,\underline X})=0$ for
every $i>0$.
\end{claim}
\begin{pfclaim} In light of theorem \ref{th_Cartan}(i) and
lemma \ref{lem_standard-coverings}, it suffices to show that
for every rational subset $R$ of $\Spa\,\underline X$, and
every standard covering $\fU$ of $R$, the alternating \v{C}ech
complex $C^\bullet_\alt(\fU,\cO^\wedge_{\Spa\,\underline X})$ is
acyclic. However, notice that any such $R$ is represented by
an affinoid (topologically local) scheme; then, arguing as
in the proof of theorem \ref{th_compl-an-noeth-are-sheaves}(i),
we are easily reduced to the case where $R=\Spa\,\underline X$.
Then, after choosing a subring of definition for $A$
we may define the scheme $V$, its formal completion
$V^\wedge$, the affine open coverings $V_\bullet$ and
$V^\wedge_\bullet$, the projection $\pi:V^\wedge\to V$ and
the morphism of schemes $\psi:X\to V$ as in the proof
of theorem \ref{th_compl-an-noeth-are-sheaves}. In light
of theorem \ref{th_coh-vanish}(i) and of the isomorphism
\eqref{eq_identify-Cech}, we are then reduced to showing
that $H^i(V^\wedge,\pi^*\psi_*\cO_{\!X})=0$ for every $i>0$.
Then, corollary \ref{cor_rompiballe}(iii) further reduces
to checking that $H^i(V,\psi_*\cO_{\!X})=0$ for every $i>0$;
but since $X$ is affine, the morphism $\psi$ is affine, so
that $H^i(V,\psi_*\cO_{\!X})\simeq H^i(X,\cO_{\!X})$, whence
the contention.
\end{pfclaim}

Now, say that $X=\Spec\,A\setminus\Spec\,A/J$ for a
finitely generated ideal $J\subset A$, and pick a finite
system of generators $f_\bullet:=(f_0,\dots,f_n)$ for
$J$ and a subring of definition $A_0\subset A$.
Let also $R_\bullet$ be the standard covering of
$\Spa\,\underline X$ associated with $f_\bullet$, and
notice that $R_i$ is represented by an affinoid scheme,
for every $i=0,\dots,n$. Set as well
$U_i:=\Spec\,A[f_i^{-1}]$ for $i=0,\dots,n$, and notice
that $\sigma_{\underline X}(R_i)\subset U_i$ for $i=0,\dots,n$,
and $X=U_0\cup\cdots\cup U_n$. There follows a commutative
diagram (notation of \eqref{subsec_altern-pseudo-Leray})
$$
\xymatrix{ A_0^\wedge\otimes_{A_0}H^i_\alt(U_\bullet,\cO_{\!X})
\ar[r]^-\alpha \ar[d] &
H^i_\alt(R_\bullet,\cO^\wedge_{\Spa\,\underline X}) \ar[d] \\
A_0^\wedge\otimes_{A_0}H^i(X,\cO_{\!X}) \ar[r] &
H^i(\Spa\,\underline X,\cO^\wedge_{\Spa\,\underline X})
}$$
where the vertical arrows are isomorphisms, by virtue of
claim \ref{cl_loud} and corollary \ref{cor_Leray}(ii). Since
$A^\wedge=A^\wedge_0\otimes_{A_0}A$, it then suffices to check
that $\alpha$ is an isomorphism. On the other hand, as
in the proof of theorem \ref{th_compl-an-noeth-are-sheaves},
we may associate with the sequence $f_\bullet$ a morphism of
schemes $\psi:X\to V$, the formal completion $V^\wedge$
of $V$ with its natural projection $\pi:V^\wedge\to V$,
and an affine covering $V_\bullet:=(V_i~|~i=0,\dots,n)$
of $V$ such that $\psi^{-1}V_i=U_i$ for $i=0,\dots,n$;
we set $V^\wedge_i:=\pi^{-1}V_i$ for $i=0,\dots,n$, and
$\cF:=\pi^*\psi_*\cO_{\!X}$. There follows a commutative
diagram
$$
\xymatrix{
A^\wedge_0\otimes_{A_0}H^i_\alt(V_\bullet,\psi_*\cO_{\!X})
\ar[r]^-\beta \ar[d] & H^i_\alt(V^\wedge_\bullet,\cF) \ar[d] \\
A^\wedge_0\otimes_{A_0}H^i(V,\psi_*\cO_{\!X}) \ar[r] &
H^i(V^\wedge,\cF)
}$$
whose vertical arrows are again isomorphisms, by
virtue of proposition \ref{prop_up-to-completion}(iv),
and the same holds for the bottom horizontal arrow, by
corollary \ref{cor_rompiballe}(iii).
To conclude, it suffices to observe that
$H^i_\alt(U_\bullet,\cO_{\!X})=H^i_\alt(V_\bullet,\psi_*\cO_{\!X})$,
and that the isomorphism of \v{Cech} complexes
\eqref{eq_identify-Cech} identifies $\alpha$ with $\beta$.
\end{proof}

\begin{example}\label{ex_counterex-sheaf-property}
We present two counterexamples that show how theorem
\ref{th_compl-an-noeth-are-sheaves} can fail for
affinoid schemes that are not analytically noetherian.
To begin with, let $K$ be a rank one complete valued
field, and fix an element $\pi\neq 0$ in the maximal
ideal of $K^+$. For every $K^+[T]$-module $J$ without
$\pi$-torsion, let us set
$$
A_J:=K^+[T]\oplus J
\qquad
A_{J,K}:=A_J\otimes_{K^+}K
$$
and endow $A_J$ with the $\pi$-adic topology and with
the $K^+[T]$-algebra structure such that the
inclusion map $K^+[T]\to A$ is a ring homomorphism
and $J$ is an ideal with $J^2=0$ ({\em i.e.}
$(x,y)\cdot(x',y'):=(xx',xy'+x'y)$ for every
$(x,y),(x',y')\in A$). Clearly $A$ is an adic topological
ring, and we endow $A_{J,K}$ with the f-adic topology such
that $A_J$ is a ring of definition. Let also $A_J^\wedge$
and $A_{J,K}^\wedge$ be the completions of $A_J$ and
respectively $A_{J,K}$, and $A^{\wedge+}_{J,K}$ the integral
closure of $A^\wedge_J$ in $A^\wedge_{J,K}$, and set
$$
\underline X:=\sSpec(A^\wedge_{J,K},A_{J,K}^{\wedge+}).
$$

$\bullet$\ \ 
Take first $J:=K^+[T,T^{-1}]/K^+[T]$. We claim that in
this case the presheaf $\cO^\wedge_{\Spa\,\underline X}$ on
$\cR(\underline X)$ is not separated. Indeed, consider the
standard covering associated with the sequence $(\pi,T)$;
it consists of the rational subsets
$$
U_0:=\{v\in\Spa\,\underline X~|~v(T)\leq v(\pi)\}
\qquad\text{and}\qquad
U_1:=\{v\in\Spa\,\underline X~|~v(\pi)\leq v(T)\}
$$
and $\cO^\wedge_{\Spa\,\underline U_0}(U_0)=B^\wedge_0\otimes_{K^+}K$,
where $B^\wedge_0$ is the completion of
$B_0:=A_J[T/\pi]\subset A_{J,K}$, where the latter is endowed
with its $\pi$-adic topology. But notice that for every
$k,n>0$ the element $T^{-k}$ of $J$ can be written in $B_0$ as
$$
T^{-k}=T^{-k-n}\cdot(T/\pi)^n\cdot\pi^n\in\pi^nB_0.
$$
Since $n$ is arbitrary, it follows that $T^{-k}$ lies
in the kernel of the restriction map
$\rho_0:A^\wedge_{J,K}\to\cO^\wedge_{\Spa\,\underline X}(U_0)$, and
since $k$ is arbitrary, we get : $J\subset\Ker\,\rho_0$.
Likewise,
$\cO^\wedge_{\Spa\,\underline U_0}(U_1)=B^\wedge_1\otimes_{K^+}K$,
where $B^\wedge_1$ is the completion of
$B_1:=A_J[\pi/T]\subset A_J[T^{-1}]$, where the latter is
endowed with its $\pi$-adic topology. Clearly
$A_J[T^{-1}]=K^+[T,T^{-1}]$, whence $J\subset
\Ker\,(\rho_1:A^\wedge_{J,K}\to\cO^\wedge_{\Spa\,\underline X}(U_1))$
as well, and the claim follows.

$\bullet$\ \ 
Next, we take $J$ to be the $\N$-graded $K^+[T]$-module
such that
$$
\gr_nJ:=T^{-n}K^+[T]/K^+[T]
\qquad
\text{for every $n\in\N$}.
$$
It follows that $\cO^\wedge_{\Spa\,\underline X}(U_0)=C_0\otimes_{K^+}K$,
where $C_0$ is the $\pi$-adic completion of $K^+[T/\pi]\oplus J'$,
with $J'$ the $\N$-graded $K^+[T/\pi]$-module such that
$$
\gr_nJ':=T^{-n}K^+[T/\pi]/K^+[T/\pi]
\qquad
\text{for every $n\in\N$}.
$$
Especially, $\cO^\wedge_{\Spa\,\underline X}(U_0)$ contains the
$\pi$-adic completion $J'^\wedge$ of $J'$. On the other hand,
it is clear that $A_J[T^{-1}]=K[T,T^{-1}]$, so the kernel of
the restriction map $\rho_{01}:\cO^\wedge_{\Spa\,\underline X}(U_0)\to
\cO^\wedge_{\Spa\,\underline X}(U_0\cap U_1)$ contains $J'^\wedge$.
Now, consider the series
$$
P:=\sum_{n\in\N}T^{-n}
\qquad
\text{where $T^{-n}\in\gr_{2n}J'$ for every $n\in\N$}.
$$
Then $T^{-n}=T^{-2n}\cdot(T/\pi)^n\cdot\pi^n\in\pi^n\gr_{2n}J'$,
so the series $P$ converges $\pi$-adically to a well defined
element of $J'^\wedge$, and $\rho_{01}(P)=0=\rho_{10}(0)$
(where $\rho_{10}:\cO^\wedge_{\Spa\,\underline X}(U_1)\to
\cO^\wedge_{\Spa\,\underline X}(U_0\cap U_1)$ is the other
restriction map). We shall show that the series $P$
does not lie in the image of $\rho_0$, so again
$\cO^\wedge_{\Spa\,\underline X}$ is not a sheaf, and more
precisely, not even the maximal separated quotient of
$\cO^\wedge_{\Spa\,\underline X}$ is a sheaf. Indeed, recall
that the completion $J^\wedge$ of $J$ lies in the product
of the completions $\gr_nJ^\wedge$ of its graded summands,
and likewise for $J'^\wedge$ (remark
\ref{rem_graded-top-algs}(iii)); moreover, the induced map
$J^\wedge\to J'^\wedge$ is the restriction of the product
of the corresponding maps $\gr_nJ^\wedge\to\gr_nJ'^\wedge$.
It is easily seen that the latter maps are injective,
so if $P$ were in the image of $\rho_0$, it would be
represented by the sequence $(P_n~|~n\in\N)$ where
$P_{2k+1}=0$ and $P_{2k}=T^{-k}\in\gr_{2k}J^\wedge$ for every
$k\in\N$. But then the series $\sum_{n\in\N}P_n$ must
converge $\pi$-adically in $J^\wedge$ (proposition
\ref{prop_Cauchy}(ii)); the latter means that for
every $t\in\N$ there exists $i\in\N$ such that
$T^{-k}\in\pi^t\gr_{2k}J$ for every $k\geq i$. But
$T^{-k}$ is not divisible by $\pi$ in $\gr_{2k}J$,
for any $k\in\N$, a contradiction.
\end{example}

\begin{remark} The constructions of example
\ref{ex_counterex-sheaf-property}, and several
others as well, may also be found in the recent preprint
\cite{Bu-Ve}, which introduces an interesting general
class of ``stably uniform'' affinoid rings $\underline A$
for which the authors can prove that the presheaf
$\cO^\wedge_{\Spa\,\underline A}$ is a sheaf.
\end{remark}

\subsection{Special loci of quasi-affinoid schemes}
The results of this section will be used in chapter
\ref{chap_perfectoid}, for the study of perfectoid spaces.

\begin{definition}\label{def_spread}
Let $\underline X:=(X,\cT_X,A^+_X)$ be any quasi-affinoid
scheme, $V\subset X^+:=\Spec\,A^+_X$ any open subset, and
set $A_X:=\cO_{\!X}(X)$. Let moreover
$\beta_{\underline X}:X\to X^+$ be the restriction of the
morphism of schemes $\gamma_{\underline X}:\Spec\,A_X\to X^+$
induced by the inclusion map $A^+_X\to A_X$.

(i)\ \
We say that $\underline X$ {\em spreads over $V$},
if $\beta_{\underline X}$ restricts to an isomorphism of
schemes $\beta^{-1}_{\underline X}V\isom V$.

(ii)\ \
Let $\underline Y:=(Y,\cT_Y,A^+_Y)$ be another quasi-affinoid
scheme, and $\phi:\underline Y\to\underline X$ an f-adic
morphism of quasi-affinoid schemes. Then the induced ring
homomorphism $\phi^\flat_Y:A_Y\to A_X$ restricts to an f-adic
map $\phi^{\flat+}_Y:A^+_Y\to A^+_X$, and we set
$$
\phi^+:=\Spec\,\phi^{\flat+}_Y:Y^+\to X^+.
$$
We say that {\em $\phi$ spreads over $V$}, if $\underline Y$
spreads over $(\phi^+)^{-1}(V)$ (in which case, we also say
that {\em $\underline Y$ spreads over $V$}, if no ambiguity
is likely to arise).

(iii)\ \
We shall also denote by $\Omega_{\underline X}\subset X^+$ the
largest open subset such that $\underline X$ spreads over
$\Omega_{\underline X}$.
\end{definition}

\begin{remark}
(i)\ \
Let $\underline X:=(X,\cT_X,A^+_X)$ be any quasi-affinoid
scheme; then $\underline X$ spreads over the analytic
locus of $\Spec\,A^+_X$, by virtue of lemma \ref{lem_deja-vu}(iii).

(ii)\ \ 
Keep the notation of definition \ref{def_spread}(ii), let
$\psi:\underline Z\to\underline Y$ be another morphism of
quasi-affinoid schemes, and suppose that $\phi$ spreads
over $V$ and $\psi$ spreads over $(\phi^+)^{-1}(V)$; then
clearly $\phi\circ\psi$ spreads over $V$.
\end{remark}

\begin{lemma}\label{lem_about-Omegas}
With the notation of definition {\em\ref{def_spread}}, we have :
\begin{enumerate}
\item
Let $W\subset X$ be an open subset such that $\beta_{\underline X}$
restricts to an open immersion $W\to X^+$. Then
$\beta_{\underline X}(W)\subset\Omega_{\underline X}$.
\item
$\gamma_{\underline X}$ restricts to an isomorphism of schemes
$\gamma_{\underline X}^{-1}\Omega_{\underline X}\isom\Omega_{\underline X}$.
\end{enumerate}
\end{lemma}
\begin{proof}(i): First, let us notice the following

\begin{claim}\label{cl_was-item-iv}
$\cO_{\!X^+}(V)$ is integrally closed in
$\cO_{\!X}(\beta_{\underline X}^{-1}V)$, for every open subset
$V\subset X^+$.
\end{claim}
\begin{pfclaim} Suppose first that $V=\Spec\,A^+_X[f^{-1}]$
for some $f\in A^+_X$; since $\beta_{\underline X*}\cO_{\!X}$ is
a quasi-coherent $\cO_{\!X^+}$-module, we have
$\cO_{\!X}(\beta^{-1}_{\underline X}V)=A_X[f^{-1}]$, and since
$A^+_X$ is integrally closed in $A_X$, the assertion
holds in this case. For a general open subset $V$, we may
find an affine open covering $V=\bigcup_{i\in I}V_i$ where
each $V_i$ is of the form $\Spec\,A^+_X[f_i^{-1}]$ for some
$f_i\in A^+_X$, and we get a commutative diagram of rings
$$
\xymatrix{ \cO_{\!X^+}(V) \ar[r] \ar[d] &
R:=\prod_{i\in I}\cO_{\!X^+}(V_i) \ar[d]
\ar@<.5ex>[r]^-{\rho_1} \ar@<-.5ex>[r]_-{\rho_2}
& \prod_{i,j\in I}\cO_{\!X^+}(V_i\cap V_j) \ar[d] \\
\cO_{\!X}(\beta_{\underline X}^{-1}V) \ar[r] &
R':=\prod_{i\in I}\cO_{\!X}(\beta^{-1}_{\underline X}V_i)
\ar@<.5ex>[r] \ar@<-.5ex>[r]
& \prod_{i,j\in I}\cO_{\!X}(\beta^{-1}_{\underline X}(V_i\cap V_j)).
}$$
Now, if $h\in\cO_{\!X}(\beta_{\underline X}^{-1}V)$ is integral over
$\cO_{\!X^+}(V)$, obviously the image $h'$ of $h$ in $R'$ is
integral over $R$; but from the foregoing case, it is easily
seen that $R$ is integrally closed in $R'$, so $h'\in R$,
and since the vertical arrows are injective maps, we deduce
easily that $\rho_1(h')=\rho_2(h')$, whence $h'\in\cO_{\!X^+}(V)$,
and the assertion follows.
\end{pfclaim}

Now, by assumption, $\beta_{\underline X}(W)$ is an open subset
of $X^+$, and $\beta_{\underline X}$ restricts to a separated
morphism
$\beta^{-1}_{\underline X}\beta_{\underline X}(W)\to\beta_{\underline X}(W)$;
the latter admits a section whose image is the open subset
$W\subset\beta^{-1}_{\underline X}\beta_{\underline X}(W)$.
Taking into account \cite[Ch.I, Cor.5.4.6]{EGAI}, we see
that $W$ is an open and closed subset of
$\beta^{-1}_{\underline X}\beta_{\underline X}(W)$.
There follow ring homomorphisms
$$
\cO_{\!X^+}(\beta_{\underline X}(W))\to
\cO_{\!X}(\beta_{\underline X}^{-1}\beta_{\underline X}(W))
\xrightarrow{\ \rho\ }\cO_{\!X}(W)
$$
whose composition is an isomorphism, and the kernel of
$\rho$ is generated by an idempotent element
$e\in\cO_{\!X}(\beta_{\underline X}^{-1}\beta_{\underline X}(W))$.
However, claim \ref{cl_was-item-iv} implies that $e$ lies
in the image of $\cO_{\!X^+}(\beta_{\underline X}(W))$, so $e=0$,
hence $W=\beta^{-1}_{\underline X}\beta_{\underline X}(W)$, and the
assertion follows immediately.

\begin{claim}\label{cl_reinstated}
Let $f:X\to Y$ and $g:Y\to Z$ be two morphisms of schemes,
such that $h:=g\circ f$ is an open immersion and the image
of $f$ is schematically dense. Then we have :
\begin{enumerate}
\item
$f$ is an open immersion.
\item
If $g$ is separated, then $g^{-1}h(X)=f(X)$.
\end{enumerate}
\end{claim}
\begin{pfclaim}(i): We know from \cite[Ch.I, Cor.5.3.13]{EGAI}
that the morphism $f$ is an immersion, and since its image
is schematically dense, it must be an open immersion.

(ii): Since $h$ is an open immersion, we may replace
$Z$ by $h(X)$ and $Y$ by $g^{-1}h(X)$, after which we
have $h=\one_X$, and $f$ is a section of $g$. Since
$g$ is separated, \cite[Ch.I, Cor.5.4.6]{EGAI} then
tells us that $f$ is a closed immersion in this case,
and combining with (i), the claim follows.
\end{pfclaim}

(ii): Consider now the commutative diagram of schemes
$$
\xymatrix{
\beta^{-1}_{\underline X}\Omega_{\underline X} \ar[r]^-f \ar[rd]_h &
\gamma_{\underline X}^{-1}\Omega_{\underline X} \ar[d]^g \\
& \Omega_{\underline X}
}$$
where $f$ is the restriction of $\beta^\circ_{\underline X}$ and
$g$ is the restriction of $\gamma_{\underline X}$. Thus, $h$ is
the restriction of $\beta_{\underline X}$, and it is an isomorphism,
by definition of $\Omega_{\underline X}$; moreover, clearly $f$
has schematically dense image. By claim \ref{cl_reinstated},
it follows that $f$ is a surjective open immersion, {\em i.e.}
it is an isomorphism, and then the same holds for $g$ as well.
\end{proof}

\begin{proposition}\label{prop_special-loci}
In the situation of definition {\em\ref{def_spread}},
suppose that $\underline X$ spreads over $V$. Then:
\begin{enumerate}
\item
The topological localization
$j_{\underline X}:\underline X{}_\loc\to\underline X$,
the topological henselization
$j'_{\underline X}:\underline X^\he\to\underline X$ and the
completion $j''_{\underline X}:\underline X^\wedge\to\underline X$
of\/ $\underline X$ spread over $V$.
\item
Let also $f_0,\dots,f_n$ be a sequence of elements of
$A_X$ that generate an ideal $I$ such that
$(\beta^{-1}_{\underline X}V)\cap\Spec\,A_X/I=\emptyset$, and set
$\underline B:=\sGamma(\underline X)(\frac{f_0,\dots,f_n}{f_0})$
(notation of \eqref{subsec_universal-property}). The induced
morphism $\underline Y:=\sSpec\,\underline B\to\underline X$
spreads over $V$.
\item
Let $\underline Y\to\underline X$ and $\underline Y'\to\underline X$
be two f-adic morphisms of quasi-affinoid schemes that spread
over $V$. The fibre product
$\underline Y\times_{\underline X}\underline Y'\to\underline X$
spreads over $V$ (see example {\em\ref{ex_fibre-prod-in-qaff-sch}}).
\item
Let $U\subset\beta_{\underline X}^{-1}V$ be a quasi-compact
open subset containing the analytic subset of $\underline X$.
The induced morphism $\phi:U\times_X\underline X\to\underline X$
spreads over $V$ (see example {\em\ref{ex_restriction-of-qaff}}).
\end{enumerate}
\end{proposition}
\begin{proof}(i): Say that $\underline X{}_\loc=(X',\cT_{X'},A^+_{X'})$,
and let $j_{\underline X}^+:X'{}^+\to X^+$ be the morphism induced by
the topological localization $A^+_X\to A^+_{X'}$.
By inspecting \eqref{subsec_localize-f-adic}, and taking into
account the natural identifications \eqref{eq_wrong-order}, we
get a cartesian diagram of schemes :
$$
\xymatrix{ X' \ar[r]^-{j_{\underline X}} \ar[d]_{\beta_{\underline X{}_\loc}}
& X \ar[d]^{\beta_{\underline X}} \\
X'{}^+ \ar[r]^-{j_{\underline X}^+} & X^+.
}$$
The assertion for $j_{\underline X}$ is an immediate consequence.
The same argument applies to $j'_{\underline X}$. Next, say that
$\underline X^\wedge=(X^\wedge,\cT_{X^\wedge},A^+_{X^\wedge})$, and
$\sGamma(\underline X)^\wedge=(A^\wedge_X,A^{\wedge+}_X,X')$.
According to corollary \ref{cor_ker-is-in-plus}(i), the
kernel $\cJ$ of the surjective map
$A^\wedge_X\to A_{X^\wedge}:=\cO_{X^\wedge}(X^\wedge)$ lies in
$A^{\wedge+}_X$, and we set $Y:=\Spec\,A^{\wedge+}_X/\cJ$; there
follows a commutative diagram of schemes
$$
\xymatrix{
& X^\wedge \ar[r] \ar[dl]_{\beta_{\underline X^\wedge}} \ar[d]^{\phi} &
X \ar[d]^{\beta_{\underline X}} \\
X^{\wedge+} \ar[r]^-\nu & Y \ar[r] & X^+
}$$
whose square subdiagram is cartesian. With this notation,
$\phi_*\cO_{\!X^\wedge}$ is a quasi-coherent $\cO_Y$-algebra,
and $A^+_{X^\wedge}=\cA(Y)$, where $\cA$ is the integral closure
of the image of $\cO_Y$ in $\phi_*\cO_{\!X^\wedge}$.
Let $V'\subset Y$ and $V''\subset X^{\wedge+}$ be the
preimages of $V$; under our assumptions, the morphism
$\phi^\flat:\cO_Y\to\phi_*\cO_{\!X^\wedge}$ restricts to an
isomorphism $\cO_{Y|V'}\isom(\phi_*\cO_{\!X^\wedge})_{|V'}$, so that
$\cA_{|V'}=(\phi_*\cO_{\!X^\wedge})_{|V'}$. The latter implies that
the morphism $\gamma_{\underline X^\wedge}$ restricts to an isomorphism
$\gamma_{\underline X^\wedge}^{-1}V''\isom V''$, and $\nu$ restricts
to an isomorphism $V''\isom V'$; lastly, since the image
of $X$ in $X^+$ contains $V$, it is clear that the image
of $X^\wedge$ in $Y$ contains $V'$, and then the image of $X^\wedge$
in $X^{\wedge+}$ contains $V''$, so $j''_{\underline X}$ spreads over $V$.

(ii): Set $B:=A_X[f_0^{-1}]$,
$B':=A^+_X[f_1/f_0,\dots,f_n/f_0]\subset B$, say that
$\underline Y=(Y,\cT_Y,B^+)$, and let $Y^+:=\Spec\,B^+$.
There follows a commutative diagram of schemes
$$
\xymatrix{ Y \ar[r]^-i \ar[d] & \Spec\,B \ar[r]^-\rho \ar[d] &
\Spec\,B' \ar[d] \\
X \ar[r] & \Spec\,A_X \ar[r] & X^+
}$$
whose left square subdiagram is cartesian.
Denote by $V'\subset\Spec\,B'$ and $V''\subset Y^+$ the
preimages of $V$, and let $g\in A^+_X$ be any element such
that $\Spec\,A^+_X[g^{-1}]\subset V$; according to lemma
\ref{lem_about-Omegas}(ii), the induced map
$A^+_X[g^{-1}]\to A_X[g^{-1}]$ is an isomorphism. Especially,
we may regard $f_0,\dots,f_n$ as elements of $C:=A^+_X[g^{-1}]$,
and by assumption $\sum_{i=0}^nf_iC=C$. It follows easily that
$f_0$ is invertible in the subring
$B'[g^{-1}]=C[f_1/f_0,\dots,f_n/f_0]$ of $B[g^{-1}]$, {\em i.e.}
$B'[g^{-1}]=B[g^{-1}]$, and since $g$ is arbitrary,
it follows that $\rho$ restricts to an isomorphism
$\rho^{-1}V'\isom V'$. Recalling now that $B^+$ is the integral
closure of $B'$ in $\cO_Y(Y)=B$, we deduce that $\gamma_{\underline Y}$
restricts as well to an isomorphism
$\gamma_{\underline Y}^{-1}V''\isom V''$. Lastly, since
$\beta_{\underline X}(X)$ contains $V$, clearly the image of
$\beta_{\underline Y}$ contains $V''$. The assertion follows
easily (details left to the reader).

(iii): Say that $\underline Y=(Y,\cT_Y,A^+_Y)$,
$\underline Y'=(Y',\cT_{Y'},A^+_{Y'})$, and set
$Y^+:=\Spec\,A^+_Y$, $Y'^+:=\Spec\,A^+_{Y'}$. Set also
$\underline Z:=\underline Y\times_{\underline X}\underline Y'$
and $(A_Z,A^+_Z,Z):=\sGamma(\underline Z)$, so that
$Z=Y\times_XY'$. Denote
$$
\phi:Z\to Z':=Y^+\times_{X^+}Y'^+
\qquad\text{and}\qquad
\pi:Z'\to X^+
$$
the induced morphisms. With this notation, $\phi_*\cO_{\!Z}$
is a sheaf of quasi-coherent $\cO_{Z'}$-algebras, and
$A^+_Z=\cA(Z')$, where $\cA$ is the integral closure
of the image of $\cO_{\!Z'}$ in $\phi_*\cO_{\!Z}$. Set
$V':=\pi^{-1}V$; under the current assumptions it is clear
that the associated morphism
$\phi^\flat:\cO_{\!Z'}\to\phi_*\cO_{\!Z}$ restricts
to an isomorphism $\cO_{\!Z'|V'}\isom(\phi_*\cO_{\!Z})_{|V'}$,
whence $\cA_{|V'}=(\phi_*\cO_{\!Z})_{|V'}$, which means that
the natural morphism $\nu:Z^+:=\Spec\,A^+_Z\to Z'$ restricts
to an isomorphism $\nu^{-1}V'\isom V'$; lastly, since
the image of $\phi$ contains $V'$, it is clear that
the image of $\beta_{\underline Z}$ contains $\nu^{-1}V'$,
whence the assertion.

(iv): Set $\underline U:=U\times_X\underline X$, and let
$h:U\to X^+$ be the composition of the open immersion
$U\to\beta^{-1}_{\underline X}V$ and the restriction
$\beta^{-1}_{\underline X}V\to X^+$ of $\beta_{\underline X}$.
Then $h=\phi^+\circ\beta_{\underline U}$, and clearly
$\beta_{\underline U}$ has schematically dense image; by
claim \ref{cl_reinstated} we deduce that $\beta_{\underline U}$
is an open immersion, whence the contention.
\end{proof}

\sset\subsubsection{}\label{subsec_spread-site}
Let $\underline X$ and $V$ be as in definition \ref{def_spread},
and suppose that $\underline X$ is topologically local and
spreads over $V$. We denote by
$$
\cR_V(\underline X)
$$
the site whose underlying category is the full subcategory
of $\cR(\underline X)$ whose objects are the rational subsets
$R$ such that the sub-presheaf $h_R\subset h_{\Spa\,\underline X}$
is represented by a morphism $\underline Y\to\underline X$
of quasi-affinoid schemes that spreads over $V$. Notice that,
in view of proposition \ref{prop_special-loci}(iii), if
$R,R'\in\Ob(\cR_V(\underline X))$, then also
$R\cap R'\in\Ob(\cR_V(\underline X))$. As usual, a family
$(R_i\to R~|~i\in I)$ of morphisms of $\cR_V(\underline X)$
generates a covering sieve of $\cR_V(\underline X)/R$ if
and only if $\bigcup_{i\in I}R_i=R$.

In light of proposition \ref{prop_special-loci}(iii), we
see that all fibre products are representable in
$\cR_V(\underline X)$, and the inclusion functor
$\iota:\cR_V(\underline X)\to\cR(\underline X)$ commutes
with fibre products. By lemma \ref{lem_crit-continuity},
it follows that $\iota$ is continuous, {\em i.e} the induced
functor on presheaves
$\iota^\wedge:\cR(\underline X)^\wedge\to\cR_V(\underline X)^\wedge$
restricts to a functor on the respective subcategories of
sheaves
$$
\tilde\iota_*:\cR(\underline X)^\sim\to\cR_V(\underline X)^\sim.
$$
We will denote as well by $\tilde\iota_*$ the corresponding
functor on abelian sheaves. For every presheaf $\cF$ on
$\cR(\underline X)$ or $\cR_V(\underline X)$, let us
denote as usual by $\cF^a$ the sheaf associated to $\cF$.

\begin{proposition}\label{prop_change-site}
In the situation of \eqref{subsec_spread-site}, let
$\cG$ be a presheaf, $\cA$ an abelian presheaf, and
$\cF$ an abelian sheaf on $\cR(\underline X)$. The
following holds :
\begin{enumerate}
\item
The natural morphism $(\iota^\wedge\cG)^a\to\iota^\wedge(\cG^a)$
is an isomorphism.
\item
The natural morphism $\iota_*\cF\to R\tilde\iota_*\cF$ is an
isomorphism in $\sD(\Z_{\cR_V(\underline X)}\Mod)$.
\item
The natural map $\check{H}^i(\cR_V(\underline X),\iota^\wedge\cA)
\to\check{H}^i(\cR(\underline X),\cA)$ is bijective for
every $i\in\N$.
\end{enumerate}
\end{proposition}
\begin{proof}(i): Define $\cF^+$ as in the proof of theorem
\ref{th_ass-sheaf}; in view of claim \ref{cl_plus-construction}
it suffices to check that the the natural morphism
$(\iota^\wedge\cG)^+\to\iota^\wedge(\cG^+)$ is an isomorphism.
For the injectivity, let $U\in\Ob(\cR_V(\underline X))$ be
any rational subset, and $s,t\in(\iota^\wedge\cG)^+(U)$ any
two sections whose images agree in $\cG^+(U)$. Then we may
find a covering family $(U_\lambda\to U~|~\lambda\in\Lambda)$
of $\cR_V(\underline X)/U$ such that $s$ and $t$ are the
classes of compatible systems $(s_\lambda~|~\lambda\in\Lambda)$,
$(t_\lambda~|~\lambda\in\Lambda)$ with
$s_\lambda,t_\lambda\in\cG(U_\lambda)$ for every $\lambda\in\Lambda$,
and the assumption means that for every $\lambda\in\Lambda$
there exists a covering family
$(U_{\lambda,i}\to U_\lambda~|~i\in I_\lambda)$ of
$\cR(\underline X)/U_\lambda$ such that
$s_{\lambda|U_{\lambda,i}}=t_{\lambda|U_{\lambda,i}}$ for every
$\lambda\in\Lambda$ and every $i\in I_\lambda$. But by
proposition \ref{prop_special-loci}(ii) and lemma
\ref{lem_standard-coverings}, the covering $U_{\lambda\bullet}$
can be refined by a covering
$(U'_{\lambda,i'}\to U_\lambda~|~i'\in I'_\lambda)$ of
$\cR_V(\underline X)/U_\lambda$, for every $\lambda\in\Lambda$.
Then $s_{\lambda|U'_{\lambda,i'}}=t_{\lambda|U'_{\lambda,i'}}$ for every
$\lambda\in\Lambda$ and every $i'\in I'_\lambda$, so that
$s=t$. Lastly, let $s\in\cG^+(U)$ be any section; by definition
we may find a covering family $(U_\lambda\to U~|~\lambda\in\Lambda)$
of $\cR(\underline X)/U$ such that $s$ is represented by a
compatible system $(s_\lambda~|~\lambda\in\Lambda)$ with
$s_\lambda\in\cG(U_\lambda)$ for every $\lambda\in\Lambda$.
But again, the covering $U_\lambda$ can be refined by a
covering of $\cR_V(\underline X)$, and we deduce that
$s$ is the image of a section of $(\iota^\wedge\cG)^+(U)$.

(ii): The assertion means that $R^i\tilde\iota_*\cF=0$ for
every $i>0$. However, recall that $R^i\tilde\iota_*\cF$ is
the sheaf associated to the presheaf that assigns to every
$U\in\Ob(\cR_V(\underline X))$ the abelian group
$H^i(\cR(\underline X)/U,\cF_{|U})$. Hence, let us fix such
a rational subset $U$, and let
$s\in H^i(\cR(\underline X)/U,\cF_{|U})$ be any element;
it suffices to exhibit a covering
$(U_\lambda\to U~|~\lambda\in\Lambda)$ of $U$ in
$\cR_V(\underline X)$ such that the image of $s$ in
$H^i(\cR(\underline X)/U_\lambda,\cF_{|U_\lambda})$
vanishes for every $\lambda\in\Lambda$. However, notice
that there exists a covering
$(U'_{\lambda'}\to U~|~\lambda'\in\Lambda')$ of $U$ in
$\cR(\underline X)$ such that the image of $s$ in
$H^i(\cR(\underline X)/U'_{\lambda'},\cF_{|U'_{\lambda'}})$
vanishes for every $\lambda'\in\Lambda'$ (this is clear,
since these cohomology groups are computed by an injective
resolution of $\cF$, which is exact in degrees $>0$).
But by proposition \ref{prop_special-loci}(ii) and
lemma \ref{lem_standard-coverings}, the covering
$U'_\bullet$ can be refined by a covering $U_\bullet$
consisting of objects of $\cR_V(\underline X)$, and
clearly such $U_\bullet$ will do.

(iii): This is an immediate consequence of the
fact -- already remarked in the foregoing -- that every
covering family of $\cR(\underline X)$ can be refined
by a covering family of $\cR_V(\underline X)$ : the details
shall be left to the reader.
\end{proof}

\sset\subsubsection{}\label{subsec_spunky}
Consider a quasi-affinoid ring $\underline A:=(A,A^+,U)$
such that $A$ is an adic (and f-adic) topological ring; let
$A_U:=\cO_{\!U}(U)$ and
$\underline U:=(U,\cT_U,A^+_U):=\sSpec\,\underline A$, and set
$$
X_A:=\Spec\,A
\qquad
X_A^+:=\Spec\,A_U^+
\qquad
X^\circ_A:=\Spec\,A_U^\circ.
$$
The induced map $A\to A_U$ is f-adic, so its image lies in
$A^\circ_U$ (lemma \ref{lem_f-adics}(iii)), and we get
a commutative diagram of schemes :
$$
\xymatrix{ U \ar[r]^-i \ar[d]_j & X_A \\
\Spec\,A_U \ar[r]^-{j'} & X_A^\circ \ar[u]_{j''}
\ar[r]^-t & X^+_A.
}$$
Let also $Z\subset X^+_A$ be the support of the $A^+_U$-module
$A^\circ_U/A^+_U$, and $Z':=X^\circ_A\setminus j'\circ j(U)$.

\begin{lemma}\label{lem_dots-on-is}
With the notation of \eqref{subsec_spunky}, the following
holds :
\begin{enumerate}
\item
$i$, $j$ and $j'\circ j$ are open
immersions.
\item
The images of $j$ and $j'$ are schematically dense.
\item
The open subset $j'\circ j(U)$ contains the analytic locus
of\/ $X^\circ_A$.
\item
$X^+_A\setminus\Omega_{\underline U}$ is the topological closure
of $Z\cup t(Z')$.
\end{enumerate}
\end{lemma}
\begin{proof}(ii): The assertion for $j'$ is clear. Next,
let $\cI$ be the kernel of the induced morphism
$j^\flat:\cO_{\Spec\,A_U}\to j_*\cO_{\!U}$; since $U$ is
quasi-compact, $j_*\cO_{\!U}$ is a quasi-coherent
$\cO_{\Spec\,A_U}$-module, so the same holds for $\cI$, and
we are reduced to checking that $\Gamma(\Spec\,A_U,\cI)=0$,
which is clear, since $j^\flat$ induces an isomorphism on
global sections.

(i): For $i$ there is nothing to prove, and the assertion
for $j$ follows from \cite[Ch.II, Prop.5.1.2]{EGAII}.
The assertion for $j'\circ j$ follows from claim
\ref{cl_reinstated}(i).

(iii): Notice that $j'$ maps the analytic locus of
$\Spec\,A_U$ isomorphically onto that of $X^\circ_A$ (lemma
\ref{lem_deja-vu}(iii)); likewise, it is easily seen that
$j''\circ j'$ restricts to an isomorphism from the analytic
locus of $\Spec\,A_U$ onto that of $X_A$. Summing up, we
conclude that $j''$ restricts to an isomorphism from the
analytic locus of $X^\circ_A$ onto that of $X_A$, whence
the assertion.

(iv): Consider the commutative diagram of schemes
$$
\xymatrix{
\beta^{-1}_{\underline U}\Omega_{\underline U} \ar[r]^-f \ar[rd]_h &
t^{-1}\Omega_{\underline U} \ar[d]^g \\
& \Omega_{\underline U}
}$$
where $f$ is the restriction of $j'\circ j$, and $g$ is the
restriction of $t$. Thus, $h$ is the restriction of
$\beta_{\underline U}$, and it is an isomorphism, by definition
of $\Omega_{\underline U}$; moreover, clearly $f$ has schematically
dense image. By claim \ref{cl_reinstated}, it follows that $f$
is a surjective open immersion, {\em i.e.} it is an isomorphism,
and then the same holds for $g$ as well. Now, if
$x\in\Omega_{\underline U}$ is any point, it follows immediately
that $x\notin Z$; moreover $g^{-1}(x)$ is the unique point of
$t^{-1}(x)$, and this point lies in $j'\circ j(U)$. This shows
that $Z\cup t(Z')\subset X^+_A\setminus\Omega_{\underline U}$.
Conversely, suppose that $x\in X^+_A$ lies neither in
$\Omega_{\underline U}$ nor the topological closure of $Z$; then
there exists an open neighborhood $V$ of $x$ in $X^+_A$ with
$V\cap Z=\emptyset$, and it follows that $t$ restricts to an
isomorphism $t^{-1}V\isom V$. In view of (i), the restriction
$\beta_{\underline U}^{-1}V\to t^{-1}V$ of $j'\circ j$
is still an open immersion, so the same holds for the
restriction $\beta_{\underline U}^{-1}V\to V$ of $\beta_{\underline U}$,
and consequently $\beta_{\underline U}(U)\cap V=
\beta_{\underline U}(\beta_{\underline U}^{-1}V)\subset\Omega_{\underline U}$,
according to lemma \ref{lem_about-Omegas}(i). We conclude that
$x\notin\beta_{\underline U}(U)$, and since $x\in t(V)$,
we finally get $x\in t(Z')$, whence (iv).
\end{proof}

\subsection{Etale coverings of quasi-affinoid schemes}
\label{sec_eta-qaff}
Let $\Sch$ be the category of schemes (in the universe $\sU$).
For every scheme $X$, let also $X_\et$ be the full subcategory
of the category $\Sch/X$, whose objects are the \'etale morphisms
$Y\to X$ with $Y$ quasi-compact and quasi-separated; notice that
every morphism in $X_\et$ is \'etale.
We endow $X_\et$ with its standard {\em \'etale topology} $J_{X,\et}$,
whose covering families are all the systems $(f_i:Y_i\to Y~|~i\in I)$
of families of morphisms of $X$-schemes such that
$\bigcup_{i\in I}f_i(Y_i)=Y$. The site $(X_\et,J_{X,\et})$ is called
the {\em \'etale site} of $X$, and is also often denoted simply by
$X_\et$, when no ambiguities are likely to arise. It is easily
seen that $X_\et$ is a lex-site if and only if $X$ is a
quasi-compact and quasi-separated scheme (the details are left
to the reader). Recall that we have a fibration $\bCov\to\Sch$
whose fibre category over every $Y\in\Ob(\Sch)$ is the category
of finite \'etale $Y$-schemes (see \eqref{sec_etale-cov}).
From the natural functor
$$
X_\et\to\Sch
\qquad
(Y\to X)\mapsto Y
$$
we deduce a fibration
$$
\bCov_{X_\et}:=X_\et\times_\Sch\bCov\to X_\et.
$$
By faithfully flat descent, it is easily seen that $\bCov_{X_\et}$
is an ind-finite stack over the site $X_\et$.

For every $k\in\N$ let $[k]:=\{0,\dots,n\}$, and set as
well $[-1]:=\emptyset$; we denote by $\sN$ the full subcategory
of the category $\Set$ with $\Ob(\sN):=\{[k]~|~k\geq-1\}$.
Then, let $\sN_{X_\et}$ be the {\em constant presheaf of categories}
on $X_\et$ with value $\sN$, {\em i.e.} the presheaf such that
$\sN_{X_\et}(Y):=\sN$ for every $Y\in\Ob(X_\et)$, and
$\sN_X(f):=\one_\sN$ for every morphism $f$ of $X_\et$.
We get an $X_\et$-cartesian functor
$$
\omega_X:\cFib(\sN_{X_\et})\to\bCov_{X_\et}
\qquad
(Y,[k])\mapsto([k]\times Y\to Y).
$$

\begin{lemma}\label{lem_standard-stuff}
The functor $\omega_X$ is $i$-covering for $i=0,1,2$.
\end{lemma}
\begin{proof} Let $f:Y\to X$ be a finite \'etale covering, and
$x\in X$; in order to check the assertion for $i=0$, we need
to exhibit an \'etale morphism $g:X'\to X$ and $k\in\N$ from
an affine scheme $X'$ such that $x\in g(X')$, with an
isomorphism of $X'$-schemes $X'\times_XY\isom[k]\times X'$.
To this aim, let $\xi:\Spec\,\kappa\to X$ be a geometric
point localized at the point $x$ (see definition
\ref{def_strict-loc}(i)); we argue by induction on the
cardinality $c$ of the set $f^{-1}(\xi)$ (see definition
\ref{def_strict-loc}(iv)). If $c=0$, we have $x\notin f(Y)$,
and then there exists an affine open neighborhood $X'$ of
$x$ in $X$ with $X'\cap f(Y)=\emptyset$; in this case we
may take for $g$ the open inclusion $X'\to X$, and $k:=0$.

Suppose next that $c>0$, and that the assertion is already
known for all schemes $X$, all $x\in X$, and all finite
\'etale coverings $f:Y\to X$ such that the cardinality
of $f^{-1}(\xi)$ is $<c$, for any geometric point $\xi$
localized at $x$. Notice that the diagonal morphism
$\Delta_{Y/X}:Y\to Y\times_XY$ is a finite \'etale morphism,
hence its image is open and closed in $Y\times_XY$. Thus,
the $Y$-scheme $Y\times_XY$ is isomorphic to the disjoint
union of $Y$ and a finite \'etale covering $h:Z\to Y$.
The geometric point $\xi$ lifts to a geometric point
$\xi_Y:\Spec\,\kappa\to Y$, and it is then easily seen
that $h^{-1}(\xi_Y)=c-1$. By inductive assumption, we find
an \'etale morphism $g_Y:Y'\to Y$ from an affine scheme $Y'$,
such that $g_Y(Y')$ contains the support of $\xi_Y$, and
such that the $Y'$-scheme $Y'\times_YZ$ is isomorphic to
$[c-2]\times Y'$. We may then take $X':=Y'$ and $g:=f\circ g_Y$.

Next, to check the assertion for $i=1$, let $k,k'\geq-1$
be two integers, $g:[k]\times X\to[k']\times X$ a
morphism of $X$-schemes, and $x\in X$; it suffices to
find an affine open neighborhood $U$ of $x$ in $X$, and
a map of sets $\phi:[k]\to[k']$ such that the restriction
$g_{|U}:[k]\times U\to[k']\times U$ of $g$ agrees with
$\phi\times U$. To this aim, it suffices to invoke
\cite[Ch.IV, Cor.17.4.7]{EGA4} : the details are left to
the reader. Lastly, it is clear from the definitions that
$\omega_X$ is a $2$-covering functor.
\end{proof}

\sset\subsubsection{}\label{subsec_do-it-better}
Denote by $\Sch^\qcqs$ the full subcategory of $\Sch$
whose objects are the quasi-compact and quasi-separated
schemes, and by $\Sch^\qcqs_\et$ the full subcategory of
$\sMorph(\Sch^\qcqs)$ whose objects are the \'etale morphisms
of schemes. The restriction
$$
t:\Sch^\qcqs_\et\to\Sch^\qcqs
$$
of the target functor $\st:\sMorph(\Sch^\qcqs)\to\Sch^\qcqs$
is a fibration, and for every $X\in\Ob(\Sch^\qcqs)$, the fibre
category $t^{-1}X$ is the previously defined category $X_\et$,
which we endow with its \'etale topology $J_{X,\et}$; we get
therefore a well defined fibred site :
$$
(\Sch^\qcqs_\et,t,J_{\bullet,\et}).
$$
For every category $I$ and every functor $F:I\to\Sch^\qcqs$ we let
$$
F_\et:=I\times_{(\Sch^\qcqs,F)}(\Sch^\qcqs_\et,t,J_{\bullet,\et})
$$
(notation of \eqref{subsec_pullback-fib-site}). Thus, $F_\et$
is a fibred lex-site over $I$, whose fibre category over every
$i\in\Ob(I)$ is naturally identified with the lex-site $(Fi)_\et$.
Moreover, we have a natural functor :
$$
F_\et\to\Sch
\qquad
(i,Y\to Fi)\mapsto Y
$$
from which we deduce a fibration over $F_\et$ :
$$
\bCov_{F_\et}:=F_\et\times_\Sch\bCov\to F_\et
$$
whose restriction to the fibre category $(Fi)_\et$ is
naturally identified with the fibration $\bCov_{(Fi)_\et}$,
for every $i\in\Ob(I)$. Explicitly, the objects of
$\bCov_{F_\et}$ are the data $(i,f:Y\to Fi,g:Z\to Y)$ where
$f$ is an \'etale morphism, $g$ is a finite \'etale morphism,
and $i\in\Ob(I)$. Let $\sN_{F_\et}$ be the constant presheaf
of categories on $F_\et$ with value $\sN$; as in
\eqref{sec_eta-qaff}, we then have a cartesian functor of
$F_\et$-fibrations :
\set\begin{equation}\label{eq_another-constant-one}
\cFib(\sN_{F_\et})\to\bCov_{F_\et}
\qquad
(i,f:Y\to Fi,[k])\mapsto(i,f:Y\to Fi,[k]\times Y\to Y).
\end{equation}
Clearly the restriction of \eqref{eq_another-constant-one}
to each fibre category $(Fi)_\et$ is naturally identified
with the functor $\omega_{Fi}$ of \eqref{sec_eta-qaff};
in view of lemma \ref{lem_standard-stuff} and proposition
\ref{prop_fibrewise-i-cov}(i) we deduce that also
\eqref{eq_another-constant-one} is $i$-covering for $i=0,1,2$.

\sset\subsubsection{}\label{subsec_return-to-qaffs}
Let $\underline X:=(X,\cT_X,A^+_X)$ be a {\em topologically
henselian\/} quasi-affinoid scheme; we attach to $\underline X$
a site $(\cQ(\underline X),J_\cQ)$ of quasi-affinoid open subsets,
as in \eqref{subsec_rational-site}.
We now define as follows a fibred lex-site over the category
$\cQ(\underline X)$. For every $U\in\Ob(\cQ(\underline X))$
choose a topologically henselian quasi-affinoid scheme
$\underline X{}^\he_U:=(X^\he_U,\cT^\he_U,A^{\he+}_U)$
representing the sub-presheaf $h'_U$ of $h'_{\underline X}$
(notation of remark \ref{rem_yoneda-rationals}(i)).
Every inclusion $U\subset U'$ of quasi-affinoid subsets of
$\Spa\,\underline X$ induces a morphism
$$
i_{UU'}:\underline X{}^\he_U\to \underline X{}^\he_{U'}
$$
of quasi-affinoid schemes (cp. \eqref{subsec_rational-site}),
and for every further inclusion $U'\subset U''$ of quasi-affinoid
open subsets of $\Spa\,\underline X$ we have
$i_{U'U''}\circ i_{UU'}=i_{UU''}$, so we get a well defined functor
$$
\underline X{}^\he_\bullet:\cQ(\underline X)\to\Sch^\qcqs
$$
whence an associated fibred site as in \eqref{subsec_do-it-better}
$$
\pi_{\underline X}:\underline X{}^\he_{\bullet,\et}:=
\cQ(\underline X)\times_{\Sch^\qcqs}(\Sch^\qcqs_\et,t,J_{\bullet,\et})\to
\cQ(\underline X).
$$
Notice that $\underline X{}^\he_{\bullet,\et}$ is a fibred lex-site,
since $\underline X{}^\he_{U,\et}$ is quasi-compact and
quasi-separated for every $U\in\Ob(\cQ(\underline X))$. We let
$(\underline X{}^\he_{\bullet,\et},J)$ be the associated total site.

$\bullet$\ \
The identity functor
$\one_{\cQ(\underline X)}:\cQ(\underline X)\to\cQ(\underline X)$
is naturally identified with the fibration
$\cFib(\sF_\bone)\to\cQ(\underline X)$, where
$\sF_\bone:\cQ(\underline X)^o\to\bCat$ is the constant
pseudo-functor with value $\bone$. The category $\bone$
admits a unique topology $\cT_\bone$, and the resulting
pair $(\bone,\cT_\bone)$ is obviously a lex-site, so we
have the fibred lex-site
$$
(\cQ(\underline X),\one_{\cQ(\underline X)},J^*_\bullet)
$$
where $J^*_U$ is the unique topology on the fibre
$\one_{\cQ(\underline X)}^{-1}(U)\isom\bone$, for every
$U\in\Ob(\cQ(\underline X))$, and we let
$(\cQ(\underline X),J^*_\cQ)$ be the resulting total site.

$\bullet$\ \
Denote by
$$
e:\cQ(\underline X)\to\underline X^\he_{\bullet,\et}.
$$
the functor that assigns to every $U\in\Ob(\cQ(\underline X))$
the final object $X^\he_U$ of the fibre category $X^\he_{U,\et}$,
and to every inclusion $U\subset U'$ the morphism $i_{UU'}$. It
is easily seen that $e$ is a $\cQ(\underline X)$-cartesian
functor. Moreover, $e$ restricts to a morphism of sites
$X^\he_{U,\et}\to(\bone,\cT_\bone)$ for every
$U\in\Ob(\cQ(\underline X))$, so $e$ is a morphism
of fibred lex-sites, and also a morphism of total sites
$$
e:(\underline X{}^\he_{\bullet,\et},J)\to(\cQ(\underline X),J^*_\cQ).
$$

$\bullet$\ \
Lastly, for every $U\in\Ob(\cQ(\underline X))$ we have an
inclusion functor
$i_U:\underline X^\he_{U,\et}\to\underline X{}^\he_{\bullet,\et}$
which identifies $\underline X^\he_{U,\et}$ with the fibre
category over $U$. This is a weak morphism of sites
$$
i_U:(\underline X{}^\he_{\bullet,\et},J)\to\underline X^\he_{U,\et}
$$
and for $U=\Spa\,\underline X$ it is even a morphism of sites
(example \ref{ex_base-with-fin-obj}). Let now $\cE$ be a stack
on the \'etale site $X_\et$; we consider the fibration :
$$
\cE_{/\cQ}:=\St(e)_*\circ\St(i_{\Spa\,\underline X})^*\cE.
$$

\begin{remark}\label{rem_fibres-of-Qslash}
(i)\ \
The fibres of the fibration $\cE_{/\cQ}$ can be described
as follows. let $c_X:\cQ(\underline X)\to\Sch^\qcqs$
be the constant functor with value $X$; we consider as well
the fibred lex-site
$$
\cQ(\underline X)\times X_\et=\cQ(\underline X)\times_{(\Sch^\qcqs,c_X)}
(\Sch^\qcqs_\et,t,J_{\bullet,\et})\to\cQ(\underline X)
$$
(see \eqref{subsec_pullback-fib-site}) whose fibre
categories are all naturally identified with $X_\et$.
Then the rule
$U\mapsto(i_{U,\Spa\,\underline X}:X^\he_U\to X)$ defines a
natural transformation $\underline X^\he_\bullet\Rightarrow c_X$
which induces, after choosing a cleavage for the fibration
$\Sch^\qcqs_\et$, a morphism of fibred lex-sites over
$\cQ(\underline X)$ :
$$
j_{\underline X}:\underline X{}^\he_{\bullet,\et}\to
\cQ(\underline X)\times X_\et
$$
(see \eqref{subsec_morph-fibsites-from-nat-tr}) which
is also a morphism on the respective total sites
$$
j_{\underline X}:(\underline X{}^\he_{\bullet,\et},J)\to
(\cQ(\underline X)\times X_\et,J').
$$
Moreover, let $l_{\underline X}:X_\et\to\cQ(\underline X)\times X_\et$
be the inclusion functor that identifies $X_\et$ with the fibre
category over the final object of $\cQ(\underline X)$. By
example \ref{ex_base-with-fin-obj}, this functor is a
morphism of sites
$$
l_{\underline X}:(\cQ(\underline X)\times X_\et,J')\to X_\et
$$
and notice that
$i_{\Spa\,\underline X}=j_{\underline X}\circ l_{\underline X}$.
Furthermore, the projection
$p_{\underline X}:\cQ(\underline X)\times X_\et\to X_\et$ is
left adjoint to $l_{\underline X}$, and from corollary
\ref{cor_lims-and-totsites}(ii) it follows easily that
$p_{\underline X}$ is a weak morphism of sites for the topology $J'$.
Combining with proposition \ref{prop_breve-for-stacks}(ii,iii),
we deduce a pseudo-natural equivalence of pseudo-functors :
$$
\St(l_{\underline X})^*\isom\St(p_{\underline X})_*.
$$
Especially, we get an equivalence of stacks :
$$
\St(i_{\Spa\,\underline X})^*\cE\isom
\St(j_{\underline X})^*\St(p_{\underline X})_*(\cE)=
\St(j_{\underline X})^*(\cQ(\underline X)\times\cE).
$$
In light of proposition \ref{prop_no-hand-waving},
we deduce a natural equivalence of stacks on
$\underline X^\he_{U,\et}$ :
$$
\St(i_U)_*\circ\St(i_{\Spa\,\underline X})^*\cE\isom
\cE_U:=\St((i_{U,\Spa\,\underline X})_\et)^*\cE
\qquad
\text{for every $U\in\Ob(\cQ(\underline X))$}.
$$
Summing up, we conclude that the fibre category over
$U$ of the fibration $\cE_{/\cQ}$ is naturally equivalent
to the fibre category $(\cE_U)_{\underline X^\he_U}$ of the
fibration $\cE_U$, over the final object $\underline X^\he_U$
of the site $\underline X^\he_{U,\et}$.

(ii)\ \
By construction, $\cE_{/\cQ}$ is a stack on the site
$(\cQ(\underline X),J^*_\cQ)$; but in fact we have :
\end{remark}

\begin{theorem}\label{th_only-for-ind-finite}
If $\cE$ is ind-finite, the fibration $\cE_{/\cQ}$ is a stack
on the site $(\cQ(\underline X),J_\cQ)$.
\end{theorem}
\begin{proof} It suffices to check that the natural functor
$$
\cE_{/\cQ}(U)\to
\Pscolim{\cS\in J_\cQ(U)^o}\sCart_{\cQ(\underline X)}(\cS,\cE_{/\cQ})
\qquad
\text{for every $U\in\Ob(\cQ(\underline X))$}
$$
is an equivalence. Since $J_\cQ(U)$ is cofiltered for the order
given by inclusion of sieves, this $2$-colimit is represented
by the colimit of the same system of categories (example
\ref{ex_filter-2-colim-in-Cat}(iv)); then arguing as in the proof
of theorem \ref{th_henselians-are-sheaves} we may assume that
$U=\Spa\,\underline X$, and replace $J_\cQ(\Spa\,\underline X)$
by its cofinal subset consisting of the sieves generated by the
standard coverings. Thus, let $f_\bullet:=(f_0,\dots,f_n)$
and $R_\bullet$ be as in remark \ref{rem_extract-proj}(i); for
every $k\in\N$ let $[k]:=\{0,\dots,k\}$, and set as well
$[-1]:=\emptyset$. We consider the category $\Sigma^+_{2,n}$
whose objects are all the maps $[j]\to[n]$ for $j=-1,0,1,2$; for
two such objects $[j]\xrightarrow{\phi}[n]\xleftarrow{\phi'}[k]$,
the morphisms $\mu:\phi\to\phi'$ are the injective non-decreasing
maps $\mu:[j]\to[k]$ such that $\phi'\circ\mu=\phi$. The unique
map $\phi_\emptyset:[-1]\to[n]$ is clearly the initial object of
$\Sigma^+_{2,n}$. Let also $\Sigma_{2,n}\subset\Sigma^+_{2,n}$ be
the full subcategory whose objects are the morphisms $[j]\to[n]$
with $j\geq 0$. Consider the functor
$$
\Phi^+:\Sigma_{2,n}^{+o}\to\cQ(\underline X)
\qquad : \qquad
([j]\xrightarrow{\phi}[n])\mapsto
R_{\phi([j])}
$$
(here $R_\Lambda\subset\Spa\,\underline X$ is defined as
in remark \ref{rem_extract-proj}(i), for every subset
$\Lambda\subset[n]$), and let
$\Phi:\Sigma_{2,n}^o\to\cQ(\underline X)$ be the restriction
of $\Phi^+$. Taking into account the discussion of
\eqref{subsec_desc-cats}, we come down to checking that the
natural functor
$$
\cE_{/\cQ}(\Spa\,\underline X)\to
\Pslim{\Sigma_{2,n}}\cE_{/\cQ}(-)\circ\Phi^o
$$
is an equivalence (notation of \eqref{subsec_we-mention}). 
Following \eqref{subsec_pullback-fib-site}, we deduce two
other fibred lex-sites
$$
\cU^\he_{f_\bullet}:=
\Sigma^{+o}_{2,n}\times_{\cQ(\underline X)}
(\underline X^\he_{\bullet,\et},J)\to\Sigma_{2,n}^{+o}
\qquad
\underline\Sigma^{+o}_{2,n}:=
\Sigma^{+o}_{2,n}\times_{\cQ(\underline X)}
(\cQ(\underline X),\one_{\cQ(\underline X)},J^*_\bullet)
$$
and according to \eqref{subsec_pullback-fib-site} and
\eqref{sec_C2-for-fibred-sites} and proposition
\ref{prop_breve-for-stacks}(ii), the projections $\Phi^+$ and
$$
\pi_\cU:\cU^\he_{f_\bullet}\to\underline X^\he_{\bullet,\et}
$$
are cocontinuous weak morphisms of sites for the topologies
of the respective total sites. Furthermore, the functor
$e\circ\Phi^+:\Sigma^{+o}_{2,n}\to X^\he_{\bullet,\et}$ factors
uniquely through a $\Sigma^{+o}_{2,n}$-cartesian functor
$$
e_\cU:\Sigma^{+o}_{2,n}\to\cU^\he_{f_\bullet}
$$
and the projection $\pi_\cU$. Just like for $e$, the functor
$e_\cU$ is a morphism of the respective total sites. Set
$\cF:=\St(\pi_\cU)_*\circ\St(i_{\Spa\,\underline X})^*\cE$;
a direct inspection yields an equivalence of stacks
$$
\cG:=\St(e_\cU)_*\cF\isom\St(\Phi^+)_*(\cE_{/\cQ})
$$
over the total site of $\underline\Sigma^{+o}_{2,n}$, and we
are reduced to checking that the natural functor
\set\begin{equation}\label{eq_Chopin}
\cG(\phi_\emptyset)\to\Pslim{\Sigma_{2,n}}\cG(-)
\end{equation}
is an equivalence. To this aim, we proceed as in
remark \ref{rem_fibres-of-Qslash}(i) : from the constant
functor $c_X:\Sigma^{+o}_{2,n}\to\Sch^\qcqs$ we deduce the
fibred site
$$
\Sigma^{+o}_{2,n}\times X_\et=\Sigma^{+o}_{2,n}\times_{(\Sch^\qcqs,c_X)}
(\Sch^\qcqs_\et,t,J_{\bullet,\et})\to\Sigma^{+o}_{2,n}
$$
whose fibre categories are naturally identified with $X_\et$,
endowed with its \'etale topology. Let
$j_{\underline X}:\underline X{}^\he_{\bullet,\et}\to
\cQ(\underline X)\times X_\et$ be as in remark
\ref{rem_fibres-of-Qslash}(i); we deduce a morphism of
fibred lex-sites
$$
j_\Sigma:=\Sigma^{+o}_{2,n}\times_{\cQ(\underline X)}j_{\underline X}:
\cU^\he_{f_\bullet}\to\Sigma^{+o}_{2,n}\times X_\et
$$
which as usual induces a morphism on total sites, and by the
same token, the projection
$$
\Sigma^{+o}_{2,n}\times_{\cQ(\underline X)}p_{\underline X}:
\Sigma^{+o}_{2,n}\times X_\et\to X_\et
$$
is a weak morphism of sites, for the topology of the total
site, and is left adjoint to the inclusion functor
$$
l_\Sigma:=\Sigma^{+o}_{2,n}\times_{\cQ(\underline X)}l_{\underline X}:
X_\et\to\Sigma^{+o}_{2,n}\times X_\et.
$$
By proposition \ref{prop_no-hand-waving} we have then an
equivalence of stacks over the total site of $\cU^\he_{f_\bullet}$ :
$$
\cF\isom\St(j_\Sigma)^*(\Sigma^{+o}_{2,n}\times\cE).
$$
Next, let $A_0\subset A:=\cO^\he_{\Spa\,\underline X}(\Spa\,\underline X)$
be a subring of definition, and $I_0\subset A_0$ an ideal of adic
definition; set $S:=\Spec\,A_0$ and $S':=\Spec\,A_0/I_0$. As in
remark \ref{rem_extract-proj}(iv,v), the sequence $f_\bullet$
induces a morphism of schemes $X\to\P^n_S$, whose schematic
image we denote by $V$, and we let $\psi:X\to V$ be the resulting
morphism of $S$-schemes; moreover, for every non-empty subset
$\Lambda\subset[n]$ we have an open subset $V_\Lambda\subset V$
which is naturally identified with
$\Spec\,A_0[f_k/f_i~|~(k,i)\in[n]\times\Lambda]$, and we set
as well $V_\emptyset:=V$. Then, for every $\Lambda\subset[n]$
let $V^\he_\Lambda$ be the henselization of $V_\Lambda$ along its
closed subscheme $V'_\Lambda:=S'\times_SV_\Lambda$; we have natural
identifications
$$
X^\he_\Lambda\isom X\times_VV^\he_\Lambda
\qquad
\text{for every $\Lambda\subset[n]$}
$$
and we denote by $\psi^\he_\Lambda:X^\he_\Lambda\to V^\he_\Lambda$
the induced projection, and by $q_\Lambda:V^\he_\Lambda\to V$ the
composition of the henselization map $V^\he_\Lambda\to V_\Lambda$
with the open immersion $V_\Lambda\to V$. We consider the functors
$$
\begin{aligned}
G&\,:\Sigma^{+o}_{2,n}\to\Sch^\qcqs
\qquad : \qquad
([j]\xrightarrow{\phi}[n])\mapsto V^\he_{\phi([j])} \\
G'&\,:\Sigma^{+o}_{2,n}\to\Sch^\qcqs
\qquad : \qquad
([j]\xrightarrow{\phi}[n])\mapsto V'_{\phi([j])}
\end{aligned}
$$
which to every morphism $\mu:\phi\to\phi'$ as in the foregoing
assign the natural morphisms
\set\begin{equation}\label{eq_op-immersion}
V^\he_{\phi'([k])}\to V^\he_{\phi([j])}
\qquad\text{and respectively}\qquad
V'_{\phi'([k])}\to V'_{\phi([j])}.
\end{equation}
As in \eqref{subsec_do-it-better} we deduce fibred lex-sites
$$
\begin{aligned}
\cV^\he_{f_\bullet}&:=
\Sigma^{+o}_{2,n}\times_{(\Sch^\qcqs,G)}(\Sch^\qcqs_\et,t,J_{\bullet,\et})
\to\Sigma^{+o}_{2,n} \\
\cV'_{f_\bullet}&:=
\Sigma^{+o}_{2,n}\times_{(\Sch^\qcqs,G')}(\Sch^\qcqs_\et,t,J_{\bullet,\et})
\to\Sigma^{+o}_{2,n}.
\end{aligned}
$$
The morphisms $(\psi^\he_\Lambda~|~\Lambda\subset[n])$,
$(q_\Lambda~|~\Lambda\subset[n])$, the closed immersions
$(V'_\Lambda\to V^\he_\Lambda~|~\Lambda\subset[n])$ and $V'\to V$,
and the open immersions $(V'_\Lambda\to V'~|~\Lambda\subset[n])$
define natural transformations
$\underline X{}^\he_\bullet\circ\Phi^+\Rightarrow G$,
$G\Rightarrow c_V$, $G'\Rightarrow G$, $c_{V'}\Rightarrow c_V$
and respectively $G'\Rightarrow c_{V'}$, where
$c_V:\Sigma^{+o}_{2,n}\to\Sch$ is the constant functor with
value $V$, and likewise for $c_{V'}$. After fixing cleavages,
these natural transformations induce morphisms of fibred
lex-sites
$$
\cU^\he_{f_\bullet}\xrightarrow{\psi_\cU}\cV^\he_{f_\bullet}
\xleftarrow{\tau_\cV}\cV'_{f_\bullet}\xrightarrow{q'}
\Sigma^{+o}_{2,n}\times V'_\et\xrightarrow{\tau_\Sigma}
\Sigma^{+o}_{2,n}\times V_\et\xleftarrow{q}\cV^\he_{f_\bullet}
$$
(see \eqref{subsec_morph-fibsites-from-nat-tr}), which are
also morphisms of sites for the topologies of the respective
total sites (proposition \ref{prop_actually-morph-of-sites}).
Define functors $e_\cV:\Sigma^{+o}_{2,n}\to\cV^\he_{f_\bullet}$ and
$e_{\cV'}:\Sigma^{+o}_{2,n}\to\cV'_{f_\bullet}$ by the rules :
$$
\phi\mapsto(\phi,V^\he_{\phi([j])})
\qquad\text{and respectively :}\qquad
\phi\mapsto(\phi,V'_{\phi[j]})
\qquad\text{for every $\phi:[j]\to[n]$}
$$
and which assign to every morphism $\mu:\phi\to\psi$ as in the
foregoing, the morphisms \eqref{eq_op-immersion}. Just
as for $e$ and $e_\cU$, both $e_\cV$ and $e_{\cV'}$ are morphisms
of the respective total sites, and we get a commutative diagram
of morphisms of sites :
\set\begin{equation}\label{eq_avot}
{\diagram V'_\et \ar[r]^-{\tau_\et} & V^\he_\et &
X_\et \ar[l]_-{\psi_\et} \\
\Sigma^{+o}_{2,n}\times V'_\et \ar[r]^-{\tau_\Sigma} \ar[u]^{l_{V'}} &
\Sigma^{+o}_{2,n}\times V_\et \ar[u]^{l_V} & \Sigma^{+o}_{2,n}\times X_\et
\ar[l]_-{\psi_\Sigma} \ar[u]_{l_\Sigma} \\
\cV'_{f_\bullet} \ar[r]^-{\tau_\cV} \ar[d]_{e_{\cV'}} \ar[u]^{q'} &
\cV^\he_{f_\bullet} \ar[d]_{e_\cV} \ar[u]^q &
\cU^\he_{f_\bullet} \ar[l]_-{\psi_\cU} \ar[d]^{e_\cU} \ar[u]_{j_\Sigma} \\
\underline\Sigma^{+o}_{2,n} \rdouble &
\underline\Sigma^{+o}_{2,n} \rdouble & \underline\Sigma^{+o}_{2,n}.
\enddiagram}
\end{equation}
where $\psi_\Sigma$ is likewise deduced from $\psi$.
Set $\cF':=\St(\psi_\cU)_*\cF$. There follows an equivalence :
$$
\cG\isom\St(e_\cV)_*\cF'.
$$
We regard the left square of the bottom row and the right square
of the central row in \eqref{eq_avot} as oriented squares of sites,
with orientations given by $\one_{e_{\cV'}}$, respectively
$\one_{q\circ\psi_\cU}$. We remark :

\begin{claim}\label{cl_from-Gabbers-work}
For every ind-finite stack $\cA$ on $\cV^\he_{f_\bullet}$,
the base change transformation
$$
\Upsilon(\St(\one_{e_{\cV'}})^\gamma_*)_\cA:
\St(e_\cV)_*(\cA)\to\St(e_{\cV'})_*\St(\tau_\cV)^*(\cA)
$$
is a natural equivalence (notation of
\eqref{subsec_base-change-map} and
\eqref{subsec_rotate-the-cube}).
\end{claim}
\begin{pfclaim} This follows by combining corollary
\ref{cor_breath-again} and \cite[Th.1' and Cor.1]{Gab}.
Notice that the proof of {\em loc.cit.} relies on
proposition 1 of the same article; for the latter,
a complete proof is available in
\cite[Exp.XX, Prop.6.3.2]{Il-Lz-Oz}. More
precisely, in the latter reference, the required results
are stated only for stacks in groupoids; however, a direct
inspection shows that the proofs work {\em verbatim} for
arbitrary ind-finite stacks.
\end{pfclaim}

Notice as well the pseudo-natural equivalences of
pseudo-functors :
\set\begin{equation}\label{eq_more-pseudo-nats-equiv}
\St(e_\cU)_*\isom\St(e_\cV)_*\circ\St(\psi_\cU)_*
\qquad
\St(\tau_\cV)^*\circ\St(q)^*\isom\St(q')^*\circ\St(\tau_\Sigma)^*.
\end{equation}

\begin{claim}\label{cl_noise-outside}
The base change transformation
$$
\Upsilon(\St(\one_{q\circ\psi_\cU})^\gamma_*):
\St(q)^*\circ\St(\psi_\Sigma)_*\to\St(\psi_\cU)_*\circ\St(j_\Sigma)^*
$$
is a pseudo-natural equivalence.
\end{claim}
\begin{pfclaim} For every subset $\Lambda\subset[n]$ we consider
the oriented square of lex-sites :
$$
\cD_\Lambda \qquad :\qquad
{\spreaddiagramcolumns{+30pt}\diagram
X^\he_{\Lambda,\et} \ar[r]^-{j_\Lambda} \ar[d]_{\psi_\Lambda}
\drtwocell\omit{_\ \ \ \ \ \ \ \ \one_{\psi\circ j'_\Lambda}} &
X_\et \ar[d]^\psi \\
V^\he_{\Lambda,\et} \ar[r]^-{q_\Lambda} & V_\et.
\enddiagram}\qquad\qquad
$$
According to corollary \ref{cor_breath-again}, it suffices
to check that the base change transformation
$$
\Upsilon(\St(\one_{\psi\circ j'_\Lambda})^\gamma_*):
\St(q_\Lambda)^*\circ\St(\psi)_*\to
\St(\psi_\Lambda)_*\circ\St(j_\Lambda)^*
$$
is a pseudo-natural equivalence for every such $\Lambda$.
However, $V^\he_\Lambda$ is the limit of a cofiltered system
of quasi-compact and quasi-separated $V$-schemes
$(V_{\Lambda,i}~|~i\in I)$ with affine transition morphisms,
and $X^\he_\Lambda$ is the limit of the induced cofiltered
system $(X_{\Lambda,i}:=X\times_VV_{\Lambda,i}~|~i\in I)$; hence
$V^\he_{\Lambda,\et}$ and $X^\he_{\Lambda,\et}$ represent the
$2$-limits of the induced systems of lex-sites
$V_{\bullet,\et}:=(V_{i,\et}~|~i\in I)$ and
$X_{\Lambda,\bullet,\et}:=(X_{\Lambda,i,\et}~|~i\in I)$. We regard
$I$ as a category (see example \ref{ex_universe}(iii)), and
we may assume that $I$ admits a final object $i_0$ such that
$V_{\Lambda,i_0}=V$; then $V_{\Lambda,\bullet,\et}$ and $X_{\Lambda,\bullet,\et}$
yield fibred lex-sites $\underline V{}_{\Lambda,\bullet,\et}$ and
$\underline X{}_{\Lambda,\bullet,\et}$ over $I$, and the system
of projections $(X_{\Lambda,i}\to V_{\Lambda,i}~|~i\in I)$ induces
a morphism of fibred lex-sites
$\underline X{}_{\Lambda,\bullet,\et}\to\underline V{}_{\Lambda,\bullet,\et}$.
Then the assertion follows from corollary \ref{cor_sMorph-is-back}.
\end{pfclaim}

Now, from proposition \ref{prop_ind-finite-and-lower-*-stacks}
we see that $\Sigma^{+o}_{2,n}\times\cE$ is an ind-finite stack,
and then the same holds for $\cF$ (theorem
\ref{th_ind-fin-and-pullback-stacks}), and also for $\cF'$,
again by proposition \ref{prop_ind-finite-and-lower-*-stacks}.
Set
$$
\cE':=\St(\tau_\et)^*\circ\St(\psi_\et)_*(\cE).
$$
From \eqref{eq_more-pseudo-nats-equiv} and claims
\ref{cl_from-Gabbers-work} and \ref{cl_noise-outside},
we deduce an equivalence of stacks :
$$
\begin{aligned}
\cG&\isom\St(e_{\cV'})_*\circ\St(q')^*\circ\St(\tau_\Sigma)^*\circ
\St(\psi_\Sigma)_*(\Sigma^{+o}_{2,n}\times\cE) \\
&\isom\St(e_{\cV'})_*\circ\St(q')^*\circ\St(l_{V'})^*(\cE').
\end{aligned}
$$
Once again, by corollary \ref{cor_lims-and-totsites}(ii) the
projection $\Sigma^{+o}_{2,n}\times V'_\et\to V'_\et$ is a weak
morphism of sites, for the topology of the total site, and is
left adjoint to $l_V$, so we get an equivalence
$$
\St(l_{V'})^*(\cE')\isom\Sigma^{+o}_{2,n}\times\cE'
$$
(proposition \ref{prop_breve-for-stacks}(ii,iii)).
Similarly, the functor $q'$ admits a left adjoint
$$
s:\cV'_{f_\bullet}\to\Sigma^{+o}_{2,n}\times V'_\et
\qquad
([j]\xrightarrow{\phi}[n],U\xrightarrow{f}V_{\phi([j])})\mapsto
(\phi,U\xrightarrow{f}V_{\phi([j])}\to V')
$$
and as explained in remark \ref{rem_never-ending-details},
the functor $s$ is cocontinuous and is a weak morphism of
sites, for the topologies of the total sites, whence
a pseudo-natural equivalence :
$$
\St(q')^*\isom\St(s)_*.
$$
Summing up, we obtain an equivalence of stacks :
\set\begin{equation}\label{eq_namo-a-casa}
\cG\isom\St(e_{\cV'})_*\circ\St(s)_*(\Sigma^{+o}_{2,n}\times\cE').
\end{equation}
Lastly, consider the functor
$$
\Psi^+:\Sigma^{+o}_{2,n}\to V'_\et
\qquad
([j]\xrightarrow{\phi}[n])\mapsto V'_{\phi([j])}
$$
and let $\Psi:\Sigma^o_{2,n}\to V'_\et$ be the restriction of
$\Psi^+$. By inspecting the definitions, we find that the
equivalence \eqref{eq_namo-a-casa} identifies the functor
\eqref{eq_Chopin} with the similar natural functor
$$
\cE'(V')\to\Pslim{\Sigma_{2,n}}\cE'(-)\circ\Psi^o.
$$
The latter is an equivalence, since $\cE'$ is a stack on
$V'_\et$, whence the theorem.
\end{proof}

\sset\subsubsection{}\label{subsec_functoriality-of-slashQ}
Let $f:\underline X'\to\underline X$ be an adic morphism of
topologically henselian quasi-affinoid schemes, and
$(U_1,\dots,U_n)$ a finite covering of $\Spa\,\underline X$
consisting of quasi-affinoid open subsets; recall that
$U'_i:=(\Spa\,f)^{-1}U_i$ is a quasi-affinoid subset of
$\Spa\,\underline U$, for every $i=1,\dots,n$ (remark
\ref{rem_yoneda-rationals}(iv)). As in
\eqref{subsec_return-to-qaffs}, define the corresponding
functors
$$
\underline X^\he_\bullet:\cQ(\underline X)\to\Sch^\qcqs
\qquad
\underline X'^\he_\bullet:\cQ(\underline X')\to\Sch^\qcqs
$$
and for every $\Lambda\subset\{1,\dots,n\}$, set
$U_\Lambda:=\bigcap_{i\in\Lambda}U_i$ (with
$U_\emptyset:=\Spa\,\underline X$) and
$U'_\Lambda:=(\Spa\,f)^{-1}U_\Lambda$. To ease notation,
let $X_\Lambda$ (resp. $X'_\Lambda$) be the scheme underlying
the quasi-affinoid scheme $\underline X^\he_{U_\Lambda}$ (resp.
$\underline X'^\he_{U'_\Lambda}$). Recall that
$\underline X^{\prime\he}_{U'_\Lambda}$ is the topological
henselization of
$\underline Z_\Lambda:=\underline X'\times_{\underline X}
\underline X^\he_{U_\Lambda}$; the composition of the
henselization morphism
$\underline X^{\prime\he}_{U'_\Lambda}\to\underline Z_\Lambda$ and
the projection $\underline Z_\Lambda\to\underline X^\he_{U_\Lambda}$
is then a morphism of quasi-affinoid schemes
$$
f_\Lambda:\underline X^{\prime\he}_{U'_\Lambda}\to
\underline X^\he_{U_\Lambda}
$$
so we have a commutative diagram of schemes :
$$
\qquad\qquad\qquad
{\diagram X'_\Lambda \ar[r]^-{f_\Lambda} \ar[d]_{p'_\Lambda} &
X_\Lambda \ar[d]^{p_\Lambda} \\
X' \ar[r]^-f & X
\enddiagram}
\qquad\text{for every $\Lambda\subset\{1,\dots,n\}$}
$$
where $p_\Lambda$ and $p'_\Lambda$ are the natural projections.
Let also $\cE$ be an ind-finite stack on $X_\et$, and set
$$
\cE':=\St(f_\et)^*\cE
\qquad
\cE_\Lambda:=\St(p_{\Lambda,\et})^*\cE
\qquad
\cE'_\Lambda:=\St(f_{\Lambda,\et})^*\cE_\Lambda
\qquad
\text{for every $\Lambda\subset\{1,\dots,n\}$}.
$$
We consider the units of adjunction :
$$
\eta:\cE\to\St(f_\et)_*\cE'
\qquad\text{and}\qquad
\eta_\Lambda:\cE_\Lambda\to\St(f_{\Lambda,\et})_*\cE'_\Lambda
\qquad\text{for every $\Lambda\subset\{1,\dots,n\}$}.
$$

\begin{lemma}\label{lem_localize-on-Q-site}
In the situation of \eqref{subsec_functoriality-of-slashQ},
suppose that the functor
$$
\eta_{\Lambda,X_\Lambda}:\cE_{\Lambda,X_\Lambda}\to\cE'_{\Lambda,X'_\Lambda}
$$
is an equivalence for every non-empty
$\Lambda\subset\{1,\dots,n\}$. Then $\eta_X:\cE_X\to\cE'_{X'}$
is an equivalence.
\end{lemma}
\begin{proof} Define the category $\Sigma^+_{2,n}$, its
full subcategory $\Sigma_{2,n}$, the fibred site
$\underline\Sigma^+_{2,n}$ and the site
$(\cQ(\underline X),J^*_\cQ)$ as in the proof of theorem
\ref{th_only-for-ind-finite}. We consider the functor
$$
\Phi^+:\Sigma_{2,n}^{+o}\to\cQ(\underline X)
\qquad : \qquad
([j]\xrightarrow{\phi}[n])\mapsto
U_{\phi([j])}
$$
and its restriction $\Phi:\Sigma_{2,n}^o\to\cQ(\underline X)$.
We attach likewise to the covering $(U'_1,\dots,U'_n)$
of $\Spa\,\underline X'$ the corresponding functors
$\Phi'^+:\Sigma_{2,n}^{+o}\to\cQ(\underline X')$ and
its restriction $\Phi'$ to $\Sigma_{2,n}^o$. Then,
from $\Psi^+:=\underline X^\he_\bullet\circ\Phi^+$ and
$\Psi'^+:=\underline X'^\he_\bullet\circ\Phi'^+$ we get
fibred sites :
$$
\cU:=\Sigma_{2,n}^{+o}\times_{(\Sch,\Psi^+)}
(\Sch^\qcqs_\et,t,J_{\bullet,\et})
\qquad
\cU':=\Sigma_{2,n}^{+o}\times_{(\Sch,\Psi'^+)}
(\Sch^\qcqs_\et,t,J_{\bullet,\et}).
$$
Next, let $e:\Sigma_{2,n}^{+o}\to\cU$ be the functor that
assigns to every $([j]\xrightarrow{\phi}[n])\in
\Ob(\Sigma_{2,n}^{+o})$ the final object $X_{\phi([j])}$ of
the site $X_{\phi([j]),\et}$; let also $i:X_\et\to\cU$ be the
inclusion functor which identifies $X_\et$ with the fibre
category of the fibration $\cU\to\Sigma^{+o}_{2,n}$, over
the final object $[-1]\in\Ob(\Sigma^{+o}_{2,n})$. Likewise
we define the functor $e':\Sigma_{2,n}^{+o}\to\cU'$ and
$i':X'_\et\to\cU'$.

With this notation, $\Phi^+$ is a weak morphism of sites
$(\cQ(\underline X),J^*_\cQ)\to\Sigma^{+o}_{2,n}$ for the
topology of the total site of $\underline\Sigma^{+o}_{2,n}$,
and likewise for $\Phi'^+$. Also, both $i$ and $e$ are
morphisms of sites for the topologies of the total sites,
and likewise for $i'$ and $e'$. Furthermore, the rule :
$([j]\xrightarrow{\phi}[n])\mapsto f_{\phi([j])}$ for
every $\phi\in\Ob(\Sigma^{+o}_{2,n})$ defines a natural
transformation
$$
\underline X'^\he_\bullet\circ\Phi^+\Rightarrow
\underline X^\he_\bullet\circ\Phi'^+
$$
which, after choosing a cleavage for the fibration
$\Sch^\qcqs_\et\to\Sch^\qcqs$, determines a morphism
of fibred sites $\phi:\cU'\to\cU$, which is again a
morphism of sites for the topologies of the respective
total sites. A direct inspection then yields a commutative
diagram of morphisms of sites :
\begin{equation}\label{eq_functor-of-Qslash}
{\diagram \Sigma^{+o}_{2,n} \ddouble &
\ar[l]_-{e'} \ar[r]^-{i_{\underline X'}} \ar[d]^g \cU' &
X'_\et \ar[d]^{f_\et} \\
\Sigma^{+o}_{2,n} & \ar[l]_-e \ar[r]^-{i_{\underline X}} \cU & X_\et.
\enddiagram}
\end{equation}
By theorem \ref{th_only-for-ind-finite}, we associate with
$\cE$ a stack $\cE_{/\cQ}$ on $(\cQ(\underline X),J_\cQ)$;
likewise, since $\cE'$ is also ind-finite (theorem
\ref{th_ind-fin-and-pullback-stacks}), we get the stack
$\cE'_{/\cQ}$ on $(\cQ(\underline X'),J'_\cQ)$. Then, arguing
as in the proof of theorem \ref{th_only-for-ind-finite} we
get equivalences of fibrations over $\Sigma^{+o}_{2,n}$ :
$$
\St(\Phi^+)_*(\cE_{/\cQ})\isom\cE_{/\Sigma}:=
\St(e)_*\circ\St(i_{\underline X})^*\cE
\qquad
\St(\Phi'^+)_*(\cE'_{/\cQ})\isom\cE'_{/\Sigma}:=
\St(e')_*\circ\St(i_{\underline X'})^*\cE'.
$$
On the other hand, from diagram \eqref{eq_functor-of-Qslash}
we deduce an equivalence of stacks
$$
\St(g)^*\circ\St(i_{\underline X})^*(\cE)\isom
\St(i_{\underline X'})^*(\cE')
$$
whence, by adjunction, a morphism of stacks :
$$
\St(i_{\underline X})^*(\cE)\to
\St(g)_*\circ\St(i_{\underline X'})^*(\cE')
$$
and after composing with $\St(e)_*$, we arrive at a
morphism of fibrations over $\Sigma^{+o}_{2,n}$ :
$$
\omega:\cE_{/\Sigma}\to\cE'_{/\Sigma}.
$$
Now, pick cleavages for $\cE_{/\Sigma}$ and $\cE'_{/\Sigma}$,
and let $\sc$ and $\sc'$ be the associated pseudo-functors;
then $\omega$ corresponds to a pseudo-natural transformation
$\omega_\bullet:\sc\Rightarrow\sc'$. Arguing as in remark
\ref{rem_fibres-of-Qslash}(i), we see that $\sc_\phi$ is
equivalent to $\cE_{\phi([j])}$ for every
$([j]\xrightarrow{\phi}[n])\in\Ob(\Sigma^{+o}_{2,n})$, and
likewise for $\sc'_\phi$, so we have for every such $\phi$
an essentially commutative diagram of categories :
$$
\xymatrix{ \sc_\phi \ar[r]^-{\omega_\phi} \ar[d] &
\sc'_\phi \ar[d] \\
\cE_{\phi([j])} \ar[r]^-{\omega^*_\phi} & \cE'_{\phi([j])}
}$$
whose vertical arrows are equivalences. Lastly, since
$\cE_{/\cQ}$ and $\cE'_{/\cQ}$ are stacks for the topologies
$J_\cQ$ and respectively $J'_\cQ$, we know that the natural
functors
$$
\cE_X\to\Pslim{\Sigma_{2,n}}\sc
\qquad
\cE'_{X'}\to\Pslim{\Sigma_{2,n}}\sc'
$$
are equivalences, hence $\omega_\emptyset:\cE_X\to\cE'_{X'}$
is identified, up to equivalences of categories and
isomorphisms of functors, with the induced functor
$$
\Pslim{\Sigma_{2,n}}\omega_\bullet:
\Pslim{\Sigma_{2,n}}\sc\to\Pslim{\Sigma_{2,n}}\sc'.
$$
To conclude the proof, it then suffices to check that
for every $\phi:[j]\to[n]$, the functor $\omega^*_\phi$
is isomorphic to the functor $\eta_{\Lambda,X_\Lambda}$, with
$\Lambda:=\phi([j])$. To this aim, we notice that the
right square subdiagram of \eqref{eq_functor-of-Qslash}
can be regarded as an oriented square of links, after
choosing adjoints and fixing adjunctions for its four
arrows; its orientation (from the bottom left corner to
the upper right corner) is given by the identity
$\one_{i_{\underline X}\circ g}:i_{\underline X}\circ g
\Rightarrow f_\et\circ i_{\underline X'}$. Then, proposition
\ref{prop_no-hand-waving} says that the base change
transformation
$$
\Upsilon(\St(\one_{i_{\underline X}\circ g})_*^\gamma):
\St(f_\et)^*\circ\St(i_{\underline X})_*\to
\St(i_{\underline X'})_*\circ\St(g)^*
$$
is a pseudo-natural equivalence. The assertion then
follows, after combining with proposition
\ref{prop_relate-adjunctions-by-bc} : the details shall
be left to the reader.
\end{proof}

\sset\subsubsection{}\label{subsec_pull-stacks-to-compl}
Resume the situation of \eqref{subsec_return-to-qaffs}, and
let also $\underline X^\wedge:=(X^\wedge,\cT^\wedge_X,A^+_{X^\wedge})$
be the completion of $\underline X$, and denote by
$\pi_X:\underline X^\wedge\to\underline X$ the {\em completion
morphism\/} of quasi-affinoid schemes (given by the counit of
the adjunction of proposition \ref{prop_complete-qaff-sch}).
The morphism $\pi_X$ induces a morphism of \'etale sites
$\pi_{X,\et}:X^\wedge_\et\to X_\et$, and for every stack $\cE$
on $X_\et$ we let
$$
\cE^\wedge:=\St(\pi_{X,\et})^*(\cE).
$$

\begin{theorem}\label{th_pass-to-completion}
With the notation of \eqref{subsec_pull-stacks-to-compl},
if $\cE$ is ind-finite, the unit of adjunction
$$
\cE_X\to\cE^\wedge_{X^\wedge}=\St(\pi_{X,\et})_*(\cE^\wedge)_X
$$
is an equivalence of categories.
\end{theorem}
\begin{proof} To begin with, let $f:\underline Y\to\underline X$
be any f-adic morphism of quasi-affinoid schemes; since
$\pi_X:\underline X^\wedge\to\underline X$ is also f-adic,
the fibre product
$\underline Y':=\underline X^\wedge\times_{\underline X}\underline Y$
is well defined (example \ref{ex_fibre-prod-in-qaff-sch}),
and we let $\underline X^\wedge\xleftarrow{p_1}\underline Y'
\xrightarrow{p_2}\underline Y$ be the natural projections.
The morphism $f^\wedge:\underline Y^\wedge\to\underline X^\wedge$
induced by $f$, and the completion morphism
$\pi_Y:\underline Y^\wedge\to\underline Y$ determine a
unique morphism of quasi-affinoid schemes
$$
g:\underline Y^\wedge\to\underline Y'
\qquad\text{such that}\qquad
p_1\circ g=f^\wedge
\quad\text{and}\quad
p_2\circ g=\pi_Y.
$$
Moreover, since $\underline Y^\wedge$ is complete, by
adjunction $g$ factors uniquely through a morphism of
quasi-affinoid schemes
$h:\underline Y^\wedge\to\underline Y'^\wedge$ and the
completion morphism $\pi_{Y'}:\underline Y'^\wedge\to\underline Y'$.

\begin{claim}\label{cl_h-is-isom}
In the foregoing situation, $h$ is an isomorphism of
quasi-affinoid schemes.
\end{claim}
\begin{pfclaim} It suffices to show that the morphisms
$p_2\circ\pi_{Y'}$ and $\pi_Y$ satisfy the same universal
property. Thus, let $\underline Z$ be a complete quasi-affinoid
schemes, and $k:\underline Z\to\underline Y$ any morphism of
quasi-affinoid schemes; the composition
$f\circ k:\underline Z\to\underline X$ factors uniquely
through a morphism $l:\underline Z\to\underline X^\wedge$
and $\pi_X$, and there exists a unique morphism
$l':\underline Z\to\underline Y'$ such that $p_1\circ l'=l$
and $p_2\circ l'=k$. Again, $l'$ factors uniquely through
a morphism $l'^\wedge:\underline Z\to\underline Y'^\wedge$
and $\pi_{Y'}$. It follows easily that $l'^\wedge$ is the
unique morphism such that $p_2\circ\pi_{Y'}\circ l'^\wedge=k$,
whence the contention.
\end{pfclaim}

\begin{claim}\label{cl_reduce-to-affinoid}
In order to prove the theorem, we may assume that
$\underline X$ is an affinoid scheme.
\end{claim}
\begin{pfclaim} Let $U_\bullet:=(U_1,\dots,U_n)$ be an open
covering of $\Spa\,\underline X$ consisting of finitely many
affinoid open subsets. (To find such $U_\bullet$, let
$A_X:=\cO_X(X)$, and $J\subset A_X$ a finitely generated ideal
such that $\Spec\,A_X/J=\Spec\,A_X\setminus X$; then we may
take a finite set of generators $f_\bullet:=(f_1,\dots,f_n)$
of $J$, and let $U_\bullet$ be the standard covering associated
with $f_\bullet$, as in \eqref{subsec_standard-cov}.)

For every $\Lambda\subset\{1,\dots,n\}$, define the
affinoid open subset $U_\Lambda\subset\Spa\,\underline X$
as in \eqref{subsec_functoriality-of-slashQ}, and notice
that $U^\wedge_\Lambda:=(\Spa\,\pi_X)^{-1}U_\Lambda$ is an
affinoid open subset of $\Spa\,\underline X^\wedge$. To
ease notation, let also
$\underline X_\Lambda:=\underline X^\he_{U_\Lambda}$ and
$\underline X'_\Lambda:=\underline X^{\wedge\he}_{U^\wedge_\Lambda}$,
for every such $\Lambda$. Let
$p_\Lambda:\underline X_\Lambda\to\underline X$ be the natural
projection, and
$\pi_\Lambda:\underline X'_\Lambda\to\underline X_\Lambda$ the
morphism of affinoid schemes deduced from $\pi_X$, as in
\eqref{subsec_functoriality-of-slashQ}, and set
$\cE_\Lambda:=\St(p_{\Lambda,\et})^*\cE$ and
$\cE'_\Lambda:=\St(\pi_{\Lambda,\et})^*\cE_\Lambda$, for every
$\Lambda\subset\{1,\dots,n\}$. By lemma
\ref{lem_localize-on-Q-site}, in order to prove the
theorem it suffices to show that the unit of adjunction
$$
\eta_\Lambda:\cE_{\Lambda,X_\Lambda}\to\cE'_{\Lambda,X'_\Lambda}
$$
is an equivalence for every non-empty
$\Lambda\subset\{1,\dots,n\}$. To this aim, let
$q_\Lambda:(\underline X_\Lambda)^\wedge\to\underline X_\Lambda$ and
$q'_\Lambda:(\underline X'_\Lambda)^\wedge\to\underline X'_\Lambda$
be the completion morphisms, and denote
by $(X_\Lambda)^\wedge$ (resp. $(X'_\Lambda)^\wedge$) the scheme
underlying $(\underline X_\Lambda)^\wedge$
($(\underline X'_\Lambda)^\wedge$); since
$(\underline X'_\Lambda)^\wedge$ is also the completion of
$\underline X^\wedge\times_{\underline X}\underline X_\Lambda$
(remark \ref{rem_loc-hens-cplete}(vi) and corollary
\ref{cor_justify}(i)), claim \ref{cl_h-is-isom} yields
an isomorphism of affinoid schemes
\set\begin{equation}\label{eq_brex-announce}
h:(\underline X_\Lambda)^\wedge\isom(\underline X'_\Lambda)^\wedge
\qquad\text{such that}\qquad
\pi_\Lambda\circ q'_\Lambda\circ h=q_\Lambda.
\end{equation}
Set $\cE^\wedge_\Lambda:=\St(q_{\Lambda,\et})^*\cE_\Lambda$ and
$\cE'^\wedge_\Lambda:=\St(q'_{\Lambda,\et})^*\cE'_\Lambda$, and denote
$$
\eta^\wedge_\Lambda:\cE_{\Lambda,X_\Lambda}\to
\cE^\wedge_{\Lambda,(X_\Lambda)^\wedge}
\qquad\text{and}\qquad
\eta'^\wedge_\Lambda:\cE'_{\Lambda,X_\Lambda}\to
\cE'^\wedge_{\Lambda,(X'_\Lambda)^\wedge}
$$
the units of adjunction. Now, up to equivalences of
categories and isomorphisms of functors, the composition
$\eta'^\wedge_\Lambda\circ\eta_\Lambda$ is naturally identified
with the counit
$$
\cE_{\Lambda,X_\Lambda}\to
(\St((\pi_\Lambda\circ q'_\Lambda)_\et)^*\cE_\Lambda)_{(X'_\Lambda)^\wedge}
$$
and the latter is an equivalence if and only if the same
holds for $\eta^\wedge_\Lambda$, due to \eqref{eq_brex-announce}.
Thus, if both $\eta^\wedge_\Lambda$ and $\eta'^\wedge_\Lambda$
are equivalences, the same will follow for $\eta_\Lambda$.
Summing up, the theorem will follow for the quasi-affinoid
scheme $\underline X$, once we have shown that the theorem
holds for the affinoid schemes $\underline X_\Lambda$ and
$\underline X'_\Lambda$, for every non-empty
$\Lambda\subset\{1,\dots,n\}$, whence the claim.
\end{pfclaim}

For every scheme $Y$ and every finite group $G$, let $G_Y$
be the group scheme $G\times Y$ over $Y$; recall that a
(right) $G$-torsor on $Y_\et$ is the datum of a surjective
\'etale morphism of schemes $Z\to Y$, together with a
morphism $\rho:Z\times_YG_Y\to Z$ that induces an action
of the group $G_Y(T)$ of $T$-sections of $G_Y$ on the
set $Z(T)$ of $T$-sections of $Z$, for every $Y$-scheme
$T$, and such that $\rho$ and the projection
$Z\times_YG_Y\to Z$ induce an isomorphism
$$
Z\times_YG_Y\isom Z\times_YZ.
$$
A morphism of $G_Y$-torsors $(Z,\rho)\to(Z',\rho')$
is a morphism of schemes $Z\to Z'$ compatible in the
obvious way with the morphisms $\rho$ and $\rho'$.
Clearly the $G$-torsors over $Y_\et$ form a category :
$$
\Tors(Y_\et,G).
$$
According to \cite[Exp.XX, Prop.6.3.2]{Il-Lz-Oz}, in
order to prove the theorem, it suffices to check that :
\begin{enumerate}
\item
for every sheaf $\cF$ on $X_\et$, the natural map
$H^0(X_\et,\cF)\!\to\!H^0(X^\wedge_\et,\tilde\pi{}_{X,\et}^*\cF)$
is bijective
\item
For every finite group $G$, and every finite morphism of
schemes $Y\to X$, the natural functor
$\Tors(Y_\et,G)\to\Tors((X^\wedge\times_XY)_\et,G)$ is an
equivalence of categories.
\end{enumerate}

Moreover, according to \cite[Exp.XII, Prop.6.5(i)]{SGA4-3},
condition (i) is equivalent to :

\begin{enumerate}
\addenu\addenu
\item
For every finite morphism of schemes $Y\to X$, the natural
map $H^0(Y_\et,\Z/2\Z)\to H^0((X^\wedge\times_XY)_\et,\Z/2\Z)$
is bijective.
\end{enumerate}

However, for every scheme $Z$, the group $H^0(Z_\et,\Z/2\Z)$
is naturally identified with the automorphism group of the
trivial $\Z/2\Z$-torsor on $Z_\et$, so (iii) follows from (ii),
for $G:=\Z/2\Z$.

Now, for every finite morphism of schemes $f:Y\to X$, let
$\underline Y:=Y\times_X\underline X$, and recall that the
morphism of quasi-affinoid schemes
$f:\underline Y\to\underline X$ is f-adic (see example
\ref{ex_finite-ext-of-affinoids}(v)), hence
$\underline Y':=\underline X^\wedge\times_{\underline X}\underline Y$
is well defined, and according to claim \ref{cl_h-is-isom} we
have an isomorphism $h:\underline Y^\wedge\isom\underline Y'^\wedge$
of quasi-affinoid schemes. Let $G$ be any finite group; we deduce
an essentially commutative diagram :
$$
\xymatrix{ \Tors(Y_\et,G) \ar[r] \ar[d] & \Tors(Y'_\et,G) \ar[d] \\
\Tors(Y^\wedge_\et,G) \ar[r] & \Tors(Y'^\wedge_\et,G)
}$$
where $Y',Y^\wedge,Y'^\wedge$ are the schemes underlying
respectively
$\underline Y',\underline Y^\wedge,\underline Y'^\wedge$.
The vertical arrows are induced by $\pi_Y:Y^\wedge\to Y$
and $\pi_{Y'}:Y'^\wedge\to Y'$, and the bottom horizontal
arrow is induced by $h$, so it is an equivalence.
Hence, condition (ii) will follow for $f:Y\to X$, if we
show that the vertical arrows of the diagram are equivalences.
Moreover, by comparing the respective universal properties,
we easily obtain a natural isomorphism of affinoid schemes
$\underline Y'\isom
(X^\wedge\times_XY)\times_{X^\wedge}\underline X^\wedge$; since
$\underline X^\wedge$ is topologically henselian, example
\ref{ex_finite-ext-of-affinoids}(vi) implies that the same
holds for $\underline Y'$. With claim
\ref{cl_reduce-to-affinoid}, we are then reduced to checking
that for every topologically henselian affinoid scheme
$(X,\cT_X,A^+_X)$ and every finite group $G$, the completion
morphism induces an equivalence
$\Tors(X_\et,G)\to\Tors(X^\wedge_\et,G)$.

Let $A_0\subset A:=\cO_X(X)$ be a subring of definition;
let also $(A_i~|~i\in I)$ be the filtered system of finite
type $A_0$-subalgebras of $A$. 
For every $i\in I$, endow $A_i$ and $A_i^+:=A_i\cap A^+_X$
with the topologies $\cT_i$ and $\cT^+_i$ induced by the
inclusions into $A$. Let $A^\wedge_0$ be the completion of
$A_0$; then the completions of $A_i$ and $A^+_i$ are
respectively $A^\wedge_i:=A^\wedge_0\otimes_{A_0}A_i$ and
$A^{\wedge+}_i:=A^\wedge_0\otimes_{A_0}A_i^+$ for every $i\in I$,
and the completion of $A$ is $A^\wedge:=A^\wedge_0\otimes_{A_0}A$,
so $A^\wedge$ is the colimit of the induced filtered
system of $A^\wedge_0$-algebras $(A^\wedge_i~|~i\in I)$.
Set as well $X_i:=\Spec\,A_i$ and $X^\wedge_i:=\Spec\,A^\wedge_i$,
and let $\cT^\wedge_i$ be the topology of $A^\wedge_i$, for every
$i\in I$; then $\underline X{}_i:=(X_i,\cT_i,A^+_i)$ is a
topologically henselian affinoid scheme (proposition
\ref{prop_quasi-affinoid}(iii)), and its completion
$\underline X{}^\wedge_i$ is $(X^\wedge_i,\cT^\wedge_i,A^{\wedge+}_i)$.

From \cite[Ch.IV, Th.8.8.2, Th.8.10.5, Th.11.2.6]{EGAIV-3}
it follows easily that the natural functors
$$
\Pscolim{i\in I}\Tors(X_{i,\et},G)\to\Tors(X_\et,G)
\qquad
\Pscolim{i\in I}\Tors(X^\wedge_{i,\et},G)\to\Tors(X^\wedge_\et,G)
$$
are equivalences of categories. Hence it suffices to check
that for every $i\in I$ the functor
$$
\Tors(X_{i,\et},G)\to\Tors(X^\wedge_{i,\et},G)
$$
is an equivalence. We may therefore assume from start that
$A$ is an $A_0$-algebra of finite type, say $A=A_0[t_1,\dots,t_n]$
for a finite sequence $t_\bullet:=(t_1,\dots,t_n)$ of elements
of $A$. For every $k=0,\dots,n$, we shall say that the sequence
$t_\bullet$ is {\em adequate up to $k$} if for every $i=1,\dots,k$,
either one of the following conditions holds :
\begin{enumerate}
\alphaenu
\item
$t_i\in A^+$
\item
$t_i\in A^\times$ and $1/t_i\in A^+$.
\end{enumerate}

\begin{claim} In order to prove the theorem, we may
assume that $t_\bullet$ is adequate up to $n$.
\end{claim}
\begin{pfclaim} Suppose that the theorem is known for every
f-adic ring $A$ generated over some topologically henselian
subring of definition $A_0$ by a sequence of elements
$(t_1,\dots,t_n)$ adequate up to $n$.
We argue by descending induction on $k\leq n$, that the
theorem then also holds for every f-adic ring generated
over any such subring $A_0$ by a sequence $(t_1,\dots,t_n)$
that is only adequate up to $k$.
If $k=n$, there is nothing to show; let then $k\leq n$
and $t_\bullet$ a sequence adequate up to $k-1$, and
suppose that we have already proved our assertion for
all sequences adequate up to $k$. Let $\underline A$ be
the affinoid ring $(A,A^+)$; we consider the affinoid rings
$$
\underline B:=\underline A\left(\frac{1,t_k}{1}\right)
\qquad
\underline C:=\underline A\left(\frac{1,t_k}{t_k}\right)
\qquad
\underline D:=\underline B\otimes_{\underline A}\underline C
$$
(notation of \eqref{subsec_universal-property} and example
\ref{ex_f-adic-push-out}(i)). Explicitly, we have
$\underline B=(B,B^+)$ and $\underline C=(C,C^+)$, where
the ring underlying the topological ring $B$ (resp. $C$)
is $A$ (resp. $A[1/t_k]$) and $B^+$ (resp. $C^+$) is the
integral closure in $B$ (resp. in $C$) of $A^+[t_k]$ (resp.
of $A^+[1/t_k]$); moreover, $B$ (resp. $C$) admits the ring
of definition $B_0:=A_0[t_k]$ (resp. $C_0:=A_0[1/t_k]$).
Recall that the topological henselization $B^\he$ of $B$
is $B^\he_0\otimes_{B_0}B$, with $B_0^\he$ the topological
henselization of $B_0$. Hence the image in $B^\he$ of
the sequence $t_\bullet$ generates the f-adic $B^\he_0$-algebra
$B^\he$, and is adequate up to $k$. Likewise we see that the
image of $t_\bullet$ generates the $C_0^\he$-algebra $C^\he$,
and then the same follows for the image of $t_\bullet$ in the
topological ring underlying $\underline D^\he$ : the details
are left to the reader. By inductive assumption, the theorem
then holds for the affinoid schemes $\sSpec\,\underline B^\he$,
$\sSpec\,\underline C^\he$ and $\sSpec\,\underline D^\he$.
But notice that $\underline B^\he$ and $\underline C^\he$
represent two affinoid open subsets $U'$ and $U''$ of
$\Spa\,\underline X$, and $\underline D^\he$ represents
$U'\cap U''$; more precisely, $(U',U'')$ is the standard
covering of $\Spa\,\underline A$ associated with the sequence
$(1,t_k)$ (see \eqref{subsec_standard-cov}); arguing as in
the proof of claim \ref{cl_reduce-to-affinoid}, we deduce
that the theorem holds for $\sSpec\,\underline A$, as required.
\end{pfclaim}

Henceforth, we assume that $t_\bullet$ is an adequate
sequence up to $n$; let $S:=\{i\leq n~|~t_i\in A^+\}$
and $T:=\{1,\dots,n\}\setminus S$, and set $t'_i:=t_i$
for every $i\in S$ and $t'_i:=1/t_i$ for every $i\in T$.
Moreover, set $t:=\prod_{i\in T}t'_i$.
The subring $A'_0:=A_0[t'_1,\dots,t'_n]$ is open and
bounded in $A$, hence it is a subring of definition of $A$,
and after replacing $A_0$ by $A'_0$ we may assume that
$A=A_0[t^{-1}]$ for an element $t\in A_0\cap A^\times$.
Let $I_0\subset A_0$ be a finite type ideal of adic
definition, and denote by $B$ the topological ring
whose underlying ring is the same as $A$, and whose
topology is the $(t,I_0)$-adic topology as defined in
\cite[Def.5.4.10(ii)]{Ga-Ra} (invoking remark
\ref{rem_I_bullet-adic}(i), it is easily seen that
the $(t,I_0)$-adic topology is a ring topology on $A$).

Since $t\in A^\times$, the ideal $t^nI_0$ is open in
$A$ for every $n\in\N$, therefore the identity
is a continuous map $A\to B$, and its completion is
a morphism of topological rings :
$$
f:A^\wedge\to B^\wedge.
$$
Let also $g:A\to A^\wedge$ and $h:B\to B^\wedge$ be the
completion maps, $I^\wedge_0$ the topological closure of
the image of $I_0$ in $A^\wedge$, and $\cT$ the
$(t,I_0^\wedge)$-adic topology on $A^\wedge$. We notice :

\begin{claim}\label{cl_new-apt}
The map $f$ is also a morphism $(A^\wedge,\cT)\to B^\wedge$
of topological rings, and induces an isomorphism of
topological rings :
$$
f^\wedge:(A^\wedge,\cT)^\wedge\isom B^\wedge.
$$
\end{claim}
\begin{pfclaim} Let $J^\wedge_0$ be the topological closure
of the image of $I_0$ in $B^\wedge$; it is clear that $f$ is
a continuous map $(A^\wedge,\cT)\to B^\wedge$, since
$f(t^nI^\wedge_0)\subset t^nJ^\wedge_0$ for every $n\in\N$.
In order to show that $f^\wedge$ is an isomorphism, we apply
the criterion of theorem \ref{th_complete-top-grps}(iii) :
first, obviously the image of $f$ is dense in $B^\wedge$,
since the image of $A$ in $B^\wedge$ is already dense.
It remains to check that $\cT$ agrees with the topology
induced from $B^\wedge$ via $f^\wedge$; we show more precisely :
$$
t^nI^\wedge_0=f^{-1}(t^nJ^\wedge_0)
\qquad
\text{for every $n\in\N$}.
$$
Indeed, since both $t^nI^\wedge_0$ and $f^{-1}(t^nJ^\wedge_0)$ are
open subgroups of $(A^\wedge,\cT)$, it suffices to check that
$g^{-1}(t^nI^\wedge_0)=g^{-1}(f^{-1}(t^nJ^\wedge_0))$; but 
$g^{-1}(t^nI^\wedge_0)=t^ng^{-1}(I^\wedge_0)=t^nI_0$, and
$g^{-1}(f^{-1}(t^nJ^\wedge_0))=h^{-1}(t^nJ^\wedge_0)=
t^nh^{-1}(J^\wedge_0)=t^nI_0$.
\end{pfclaim}

Since $h=f\circ g$, it suffices to check that $f$ and $h$
induce equivalences of categories
$$
\Tors((\Spec\,A^\wedge)_\et,G)\isom\Tors((\Spec\,B^\wedge)_\et,G)
\xleftarrow{\sim}\Tors((\Spec\,A)_\et,G).
$$
Now, the pair $(A_0,I_0)$ is henselian by assumption, and
since $A_0^\wedge$ is complete for the $I_0^\wedge$-adic
topology, also the pair $(A^\wedge_0,I^\wedge_0)$ is henselian;
it then follows that the pairs $(A_0,tI_0)$ and
$(A^\wedge_0,tI^\wedge_0)$ are henselian
\cite[Rem.5.1.10(iv)]{Ga-Ra}. By \cite[Prop.5.4.53]{Ga-Ra}
and claim \ref{cl_new-apt}, we deduce that $f$ and $h$ induce
equivalences of categories of finite \'etale morphisms
\set\begin{equation}\label{eq_poor-myriam}
\bCov(\Spec\,A^\wedge)\isom\bCov(\Spec\,B^\wedge)\xleftarrow{\sim}
\bCov(\Spec\,A)
\end{equation}
(notation of \eqref{sec_etale-cov}). Now, for any scheme
$Y$, the $G$-torsors over $Y_\et$ are the objects $\phi:Z\to Y$
of $\bCov(\Spec\,A^\wedge)$ endowed with a (right) $G$-action
and such that $\phi$ is a surjective morphism. Such a morphism
has an open and closed image in $Y$, hence it is surjective
if and only if it is an epimorphism (in the category $\bCov(Y)$).
The equivalences \eqref{eq_poor-myriam} induce formally
equivalences for the corresponding categories of objects
endowed with $G$-action (details left to the reader), and
obviously respect epimorphisms, whence the contention.
\end{proof}

\begin{remark} For every finite group $G$, let $\cC_G$
be the category with $\Ob(\cC_G)=\{o\}$ and such that
the monoid $\Hom_{\cC_G}(o,o)$ is isomorphic to $G$.
For every scheme $X$, let also $\cC_{G,X_\et}$ be the
constant presheaf on $X_\et$ with value $\cC_G$. By
arguing as in the proof of lemma \ref{lem_standard-stuff},
it is easily seen that the associated stack $\cC_{G,X_\et}^a$
is naturally equivalent to the stack of $G$-torsors
on $X_\et$ : the details shall be left to the reader.
\end{remark}

\sset\subsubsection{}
For every $U\in\Ob(\cQ(\underline X))$, choose a complete and
separated quasi-affinoid scheme
$\underline X^\wedge_U:=(X^\wedge_U,\cT^\wedge_U,A^{\wedge+}_U)$
representing the sub-presheaf
$h''_U\subset h''_{\underline X^\wedge}$ (notation of remark
\ref{rem_yoneda-rationals}(i)). Any inclusion $U\subset U'$
of quasi-affinoid subsets of $\Spa\,\underline X$ induces a
morphism $j_{UU'}:\underline X^\wedge_U\to\underline X^\wedge_{U'}$
of quasi-affinoid schemes (cp. \eqref{subsec_rational-site}),
and for any further inclusion $U'\subset U''$ we
have
$$
j_{U'U''}\circ j_{UU'}=j_{UU''}
$$
whence a well defined functor
$$
\underline X{}^\wedge_\bullet:\cQ(\underline X)\to\Sch
\qquad
U\mapsto X{}^\wedge_U
$$
and we consider the fibration over $\cQ(\underline X)$
$$
\bCov^\wedge_{\underline X}:=\Fib(\underline X{}^\wedge_\bullet)^*(\bCov)
$$
where the fibration $\bCov\to\Sch$ is defined as in
\eqref{sec_eta-qaff}.  With this notation, we can now state :

\begin{corollary}\label{cor_Cov-is-a-stack}
The fibration $\bCov^\wedge_{\underline X}$ is a stack on the site
$(\cQ(\underline X),J_\cQ)$.
\end{corollary}
\begin{proof} We have already remarked that the stack
$\bCov_{X_\et}$ is ind-finite on $X_\et$, hence the fibration
$(\bCov_{X_\et})_{/\cQ}$ is a stack on $(\cQ(\underline X),J_\cQ)$,
by theorem \ref{th_only-for-ind-finite}. But from theorem
\ref{th_pass-to-completion} we easily deduce an equivalence
of fibrations $(\bCov_{X_\et})_{/\cQ}\isom\bCov^\wedge_{\underline X}$,
whence the assertion.
\end{proof}

\begin{remark}
A variant of corollary \ref{cor_Cov-is-a-stack} appears
in \cite[Th.2.6.9]{Ked-Liu}.
\end{remark}

\sset\subsubsection{}\label{subsec_Florensky}
We wish to give an explicit description, up to natural
equivalence of categories, of the stalks of the fibration
$\bCov^\wedge_{\underline X}$. Hence, let $x\in\Spa\,\underline X$
be any point; recall that the stalk of $\bCov^\wedge_{\underline X}$
over $x$ is the category
$$
\bCov^\wedge_{\underline X}(x):=
\Pscolim{U\in\cQ(\underline X,x)^o}\bCov(X^\wedge_U)
$$
where $\cQ(\underline X,x)$ denotes the set of all quasi-affinoid
open neighborhoods of $x$ in $\Spa\,\underline X$, which is
cofiltered by inclusion. Pick a finitely generated ideal
$J\subset A_X:=\cO_X(X)$ such that $X=\Spec\,A_X\setminus\Spec\,A_X/J$,
and let $f_0,\dots,f_n$ be a finite system of generators of $J$.
Set $U_i:= R_{A_X}\bigl(\frac{f_0,\dots,f_n}{f_i}\bigr)\cap
\Spa\,\underline X$ for $i=0,\dots,n$. Then we have
$$
U_0\cup\cdots\cup U_n=\Spa\,\underline X
$$
and we may assume that $x\in U_0$. Notice that $h''_{U_0}$
is represented by a quasi-affinoid scheme
$\underline X^\wedge_{U_0}:=
(X^\wedge_{U_0},\cT^\wedge_{U_0},A^{\wedge+}_{U_0})$ with
$X^\wedge_{U_0}=\Spec\,\cO^\wedge_{\Spa\,\underline X}(U_0)$. Then,
since $X^\wedge_{U_0}$ is affine, every rational subset
of $\Spa\,\underline X^\wedge_{U_0}$ containing $x$ is
likewise represented by a quasi-affinoid scheme whose
underlying scheme is affine, and the system of such
rational subsets is naturally identified with a cofinal
subset of $\cQ(\underline X,x)$. Taking into account
lemma \ref{lem_go-to-pseudo-lim}, we deduce a natural
equivalence
$$
\bCov^\wedge_{\underline X}(x)\isom
\bCov(\Spec\,\cO^\wedge_{\Spa\,\underline X,x}).
$$

$\bullet$\ \ 
Next, if $x$ is analytic, combining with lemma
\ref{lem_stalks-are-local}(i,v.c), proposition
\ref{prop_krasner}(iii) and \cite[Ch.IV, Prop.18.5.15]{EGA4}
we arrive at a natural equivalence
$$
\bCov^\wedge_{\underline X}(x)\isom
\bCov(\Spec\,\kappa(x^\wedge))\isom
\bCov(\Spec\,\kappa(x)^\wedge)
$$
where $\kappa(x^\wedge)$ is the residue field of the
natural valuation $|\cdot|^\wedge_x$ on
$\cO^\wedge_{\Spa\,\underline X,x}$, and $\kappa(x)^\wedge$ is the
completion of $\kappa(x^\wedge)$ for its valuation topology
(see \eqref{subsec_valuation-on-stalks}).

$\bullet$\ \
Lastly, if $x$ is non-analytic, lemma
\ref{lem_stalks-are-local}(vi.b) and \cite[Th.5.5.7(iii)]{Ga-Ra}
yield the natural equivalence
$$
\bCov^\wedge_{\underline X}(x)\isom
\bCov(\Spec\,\cO^\wedge_{\Spa\,\underline X,x}/
A^{\circ\circ}_X\cO^\wedge_{\Spa\,\underline X,x}).
$$

\section{Perfectoid rings and perfectoid spaces}
\label{chap_perfectoid}
In this chapter we develop a generalization of Scholze's
theory of perfectoid rings and perfectoid spaces.

\subsection{Distinguished elements and transversal pairs}
\label{sec_def-A-tilde}
Let $p$ be a prime number, $(A,\cT)$ a complete and
separated topological ring whose topology $\cT$ is
linear and coarser than the $p$-adic topology.
Remark \ref{rem_fontaine}(ii) endows the topological
monoid $\bE(A)$ with a natural structure of perfect
topological $\F_p$-algebra, such that :
\begin{itemize}
\item
For every topologically nilpotent ideal $I\subset A$
containing $pA$, and which is either closed for the
topology $\cT$, or else finitely generated and closed
for the $p$-adic topology of $A$, the projection
$\pi_I:A\to A/I$ induces an isomorphism of topological rings
$$
\bE(\pi_I):\bE(A)\isom\bE(A/I)
$$
where $A/I$ is endowed with the quotient topology induced
by $\cT$ via $\pi_I$.
\item
Moreover, if $(A',\cT')$ is any other topological ring
fulfilling the above conditions, and $f:A\to A'$ is any
continuous ring homomorphism, the map $\bE(f):\bE(A)\to\bE(A')$
is a morphism of topological rings.
\end{itemize}
In light of this canonical identification of $\bE(A)$ with
$\bE(A/I)$, we shall henceforth write slightly abusively
$$
\bar u_{A/I}:\bE(A)\to A/I
$$
for the map that would be denoted
$\bar u_{A/I}\circ\bE(\pi_I)$, or equivalently,
$\pi_I\circ\bar u_A$, with the notation of section
\ref{sec_Fontaine-only}. Another basic construction
of this chapter shall be the topological ring :
$$
\bA(A,\cT):=W(\bE(A))
$$
where $W(-)$ denotes the ring of Witt vectors as in
\eqref{subsec_Witt-vectors}. From remark
\ref{rem_topology-of-E}(i) and lemma
\ref{lem_Witt-limit}(ii,iii) we also know that the
topologies of $\bE(A)$ and $\bA(A,\cT)$ are linear,
complete and separated. Unless we have to deal
with several different topologies on $A$, we shall
usually omit mentioning explicitly $\cT$, and write
simply $\bA(A)$. According to proposition
\ref{prop_Witt-is-complete}(ii), the ring $\bA(A)$ is
also complete and separated for its $p$-adic topology,
and the same holds for $A$, by lemma \ref{lem_fontaine};
moreover we see from \eqref{eq_F-V-commute} that $p$
is a regular element of $\bA(A)$, and by remark
\ref{rem_Witt-limit}(i) the ghost component $\bomega_0$
descends to an isomorphism of topological rings
$$
\bA(A)\otimes_\Z\F_p\isom\bE(A)
$$
provided we endow $\bA(A)\otimes_\Z\F_p$ with the quotient
topology induced from $\bA(A)$ via the projection
$\bA(A)\to\bA(A)\otimes_\Z\F_p$. The map $\bomega_0$ also
admits a continuous set-theoretic multiplicative splitting,
the Teichm{\"u}ller mapping of \eqref{subsec_Teich}, which
will be here denoted :
$$
\tau_A:\bE(A)\to\bA(A).
$$ 

\begin{lemma}\label{lem_drop-conditions}
With the notation of \eqref{sec_def-A-tilde}, the
following holds :
\begin{enumerate}
\item
There exists a unique ring homomorphism $u_A:\bA(A)\to A$
that makes commute the diagram :
$$
\xymatrix@C+20pt{
\bA(A) \ar[r]^-{u_A} \ar[d]_{\bomega_0} & A \ar[d]^{\pi_{pA}} \\
\bE(A) \ar[r]^-{\bar u_{A/pA}} & A/pA.
}$$
\item
The map $u_A$ is continuous for the topology $\cT$ and
the topology of $\bA(A,\cT_A)$.
\item
$\bA(A)^\times=u_A^{-1}(A^\times)$.
\item
$u_A\circ\tau_A=\bar u_A$.
\end{enumerate}
\end{lemma}
\begin{proof} We have already remarked that $A$ is
complete and separated for its $p$-adic topology, and
the filtration $(J_n:=p^nA~|~n\in\N)$ of $A$ trivially
fulfills the condition of lemma \ref{lem_basic-cong};
therefore (i) follows immediately from proposition
\ref{prop_lift-Witt}(iii). More precisely, \eqref{eq_explicit}
translates as the following explicit expression for $u_A$ :
$$
u(\underline a)=\sum_{n\in\N}p^n\cdot\bar u_A(a_n^{p^{-n}})
\qquad
\text{for every $\underline a:=(a_n~|~n\in\N)\in\bA(A)$}
$$
from which (ii) follows easily : details left to the
reader.

(iii): The inclusion $\bA(A)^\times\subset u_A^{-1}(A^\times)$
is obvious. For the converse, let $\alpha\in\bA(A)$ be any
element such that $u_A(\alpha)\in A^\times$; then
$\bomega_0(\alpha)\in\bar u_A^{\ -1}(A^\times)$.
By remark \ref{rem_general}(iii) we deduce
that $\bomega_0(\alpha)\in\bE(A)^\times$, so it
remains only to notice that
$\bomega_0^{-1}(\bE(A)^\times)=\bA(A)^\times$,
due to proposition \ref{prop_Witt-is-complete}(iii).

(iv) is clear from the foregoing explicit formula for $u_A$.
\end{proof}

Furthermore, $\bA(A)$ is endowed with an automorphism
$\bsigma_{\!A}:\bA(A)\to\bA(A)$ that lifts the Frobenius map of
$\bE(A)$, {\em i.e.} such that the horizontal arrows of the
diagram :
\set\begin{equation}\label{eq_up-and-down}
{\diagram
\bA(A) \ar[rr]^-{\bsigma_{\!A}} \ar@<-.5ex>[d]_{\bomega_0} & &
\bA(A) \ar@<-.5ex>[d]_{\bomega_0} \\
\bE(A) \ar[rr]^-{\Phi_{\bE(A)}} \ar@<-.5ex>[u]_{\tau_A}
& & \bE(A) \ar@<-.5ex>[u]_{\tau_A}
\enddiagram}\end{equation}
commute with the downward arrows (proposition
\ref{prop_V_A-and-F_A}(iii)). The horizontal arrows commute
also with the upward ones, due to \eqref{eq_Teich-Frob}.

\sset\subsubsection{}\label{subsec_new-place}
Keep the notation of \eqref{sec_def-A-tilde}, and
recall that $\tau_A$ induces a continuous map
$\Spec\,\tau_A:\Spec\,\bA(A)\to\Spec\,\bE(A)$ as in
proposition \ref{prop_like-a-ring-map}(ii); in light
of lemma \ref{lem_drop-conditions}(iv) and  remark
\ref{rem_strongly-spectral}(ii) we deduce a continuous
and spectral map
$$
\Spec\,\bar u_A:=(\Spec\,\tau_A)\circ(\Spec\,u_A):
\Spec\,A\to\Spec\,\bE(A)
\qquad
\fp\mapsto\bar u{}^{\ -1}_A\fp.
$$
Moreover, due to proposition \ref{prop_like-a-ring-map}(iv),
every morphism $\phi:A\to A'$ of topological algebras that
fulfill the conditions of \eqref{sec_def-A-tilde}
induces a commutative diagram of topological spaces :
\set\begin{equation}\label{eq_functoriality-of-Spec_u}
{\diagram \Spec\,A' \ar[rr]^-{\Spec\,\bar u_{A'}}
\ar[d]_{\Spec\,\phi} & &
\Spec\,\bE(A') \ar[d]^{\Spec\,\bE(\phi)} \\
\Spec\,A \ar[rr]^-{\Spec\,\bar u_A} & & \Spec\,\bE(A).
\enddiagram}
\end{equation}

\begin{lemma}\label{lem_new-place}
In the situation of \eqref{subsec_new-place},
let $\fp\in\Spec\,A$ be any prime ideal, set
$\fq:=\Spec\,\bar u_A(\fp)\in\Spec\,\bE(A)$, and let
$\pi_\fp:A\to\kappa(\fp)$, $\pi_\fq:\bE(A)\to\kappa(\fq)$
be the projections. If\/ $\fp$ is a closed subset in the
$p$-adic topology of\/ $A$, then $\bar u_A$ induces a
morphism of multiplicative monoids
$$
\bar u_\fp:\kappa(\fq)\to\kappa(\fp)
\qquad\text{such that}\qquad
\bar u_\fp\circ\pi_\fq=\pi_\fp\circ\bar u_A.
$$
\end{lemma}
\begin{proof} Since $u_A$ is continuous for the $p$-adic
topologies, $u_A^{-1}\fp$ is a closed subet for the $p$-adic
topology of $\bA(A)$, so the assertion follows immediately
from proposition \ref{prop_like-a-ring-map}(iii).
\end{proof}

With the following definition, we single out certain
elements of $\bA(A)$ that shall play a special role in
the study of perfectoid rings.

\begin{definition}\label{def_distinguished}
Let $(E,\cT)$ be any topological $\F_p$-algebra.
\begin{enumerate}
\item
An element $(a_n~|~n\in\N)\in W(E)$ is called
{\em distinguished} if $a_0$ is topologically
nilpotent in $E$, and $a_1\in E^\times$.
\item
An ideal of $W(E)$ is called {\em distinguished}
if it is generated by a distinguished element.
\end{enumerate}
\end{definition}

\begin{remark}\label{rem_distinguished}
Let $E$ be any topological $\F_p$-algebra,
$\underline a:=(a_n~|~n\in\N)$ and
$\underline b:=(b_n~|~n\in\N)$ two elements of $W(E)$,
such that $\underline a$ is distinguished and set
$\underline c:=\underline a\cdot\underline b$. 

(i)\ \
Suppose that the topology of $E$ is linear, complete
and separated. Then $\underline c$ is distinguished if
and only if $\underline b$ is invertible in $W(E)$.
Indeed, say that $\underline c=(c_n~|~n\in\N)$. Then
remark \ref{rem_homogeneous-laws} says that $c_0=a_0b_0$,
especially $c_0$ is topologically nilpotent. Also,
$c_1=a_0^pb_1+a_1b_0^p$, and notice that $b_0\in E^\times$,
if and only if $\underline b\in W(E)^\times$ (proposition
\ref{prop_Witt-is-complete}(iii)); on the other hand, since
the topology of $E$ is linear, $a_0^pb_1$ is topologically
nilpotent (remark \ref{rem_something-on-bdd}(v)). Since
the topology of $E$ is linear, complete and separated, it
follows that $c_1$ is invertible if and only if the same
holds for $b_0$ (details left to the reader), whence the
assertion.

(ii)\ \
If $E$ is reduced, then every distinguished ideal of
$W(E)$ is a closed subset for the $p$-adic topology
of $W(E)$ and a free $W(E)$-module of rank one
(proposition \ref{prop_reduced-Witt}(i,v)).

(iii)\ \
Let $E$ be as in (i), and $I,J\subset W(E)$ any two
distinguished ideals such that $I\subset J$; then (i)
implies that $I=J$. By the same token, we see that
every generator of $I$ is distinguished. Indeed, let
$\underline a_1$ and $\underline a_2$ be any two
generators of $I$; then there exists
$\underline u:=(u_n~|~n\in\N)\in W(E)$ such that
$\underline a_2=\underline u\cdot\underline a_1$.
Let $N\subset E$ be the nilpotent ideal, and endow
$E_\red:=E/N$ with the quotient topology induced from
$E$; it follows from (ii) that the image of
$\underline u$ is invertible in $W(E_\red)$, therefore
the image of $u_0$ is invertible in $E_\red$, so $u_0$
is invertible in $E$, and finally $\underline u$ is
invertible in $W(E)$ (proposition
\ref{prop_Witt-is-complete}(iii)).
\end{remark}

\begin{example}
With the notation of example \ref{ex_Witt-ring},
notice that the only distinguished elements of
$W(A,\cT_T)=\Z_p\{T^{1/p^\infty}\}$ are of the form $pu$,
where $u\in W(A,\cT_T)^\times$ is an arbitrary element;
especially, the only distinguished ideal of
$\Z_p\{T^{1/p^\infty}\}$ is the principal ideal generated
by $p$. On the other hand, the distinguished elements
of $W(A^\wedge,\cT^\wedge_T)$ are all those of the form
$pu+aT^\lambda$, where $a$ (resp. $u$) is an arbitrary
element (resp. invertible element) of
$W(A^\wedge,\cT^\wedge_T)$, and $\lambda\in\N[1/p]$ is
any strictly positive number.
\end{example}

The following result shall be applied to produce useful
distinguished elements in various situations.

\begin{lemma}\label{lem_before-name}
Let $E$ be an $\F_p$-algebra, $\cI\subset W(E)$ an ideal,
$\underline\alpha:=(\alpha_n~|~n\in\N)\in\cI$ any element,
$\beta_\bullet:=(\beta_1,\dots,\beta_k)$ a finite system of
elements of $E$. Denote by $J$ the ideal generated by
$\beta_\bullet$, and by $u:W(E)\to A:=W(E)/\cI$ the natural
projection. Suppose furthermore that :
\begin{enumerate}
\alphaenu
\item
$\alpha_0$ lies in the Jacobson radical of $E$ and
$\alpha_1\in E^\times$.
\item
The Frobenius endomorphism $\Phi_E$ of $E$ is surjective.
\end{enumerate}

{\em (i)}\ \ Then the following conditions are equivalent :
\begin{enumerate}
\alphaenu
\addenu\addenu
\item
$pA$ is contained in the ideal of $A$ generated by
$(u\circ\tau_E(\beta_1),\dots,u\circ\tau_E(\beta_k))$.
\item
There exists an element
$\underline\alpha':=(\alpha'_n~|~n\in\N)\in\cI$ such
that $\alpha'_0\in J$ and the image of $\alpha'_1$
is invertible in $E/J$.
\item
There exists an element
$\underline\alpha':=(\alpha'_n~|~n\in\N)\in\cI$ such
that $\alpha'_0\in J$ and $\alpha'_1\in E^\times$.
\end{enumerate}

{\em(ii)}\ \ 
Suppose furthermore that $\cI=\underline\alpha W(E)$.
Then {\em(c),(d)} and {\em(e)} are also equivalent to :
\begin{enumerate}
\alphaenu
\addenu\addenu\addenu\addenu\addenu
\item
$\alpha_0\in J$.
\end{enumerate}
\end{lemma}
\begin{proof}(c)$\Rightarrow$(d): Assumption (c) implies
that there exist
$\underline\gamma{}_1,\dots,\underline\gamma{}_k\in W(E)$ such
that
$$
\underline\alpha':=
p+\sum_{i=1}^k\underline\gamma{}_i\cdot\tau_E(\beta_i)\in\cI
$$
and we claim that (d) holds for this $\underline\alpha'$.
Indeed, $\alpha'_0=
\sum_{i=1}^k\bomega_0(\underline\gamma{}_i)\cdot\beta_i\in J$,
and clearly the image of $\underline\alpha'$ in $W(E/J)$
equals $p$, {\em i.e.} $\underline\alpha'_1\equiv 1\pmod{J}$.

(d)$\Rightarrow$(c): Let $\underline\alpha'$ be as in (d);
then there exists $\lambda\in E$ such that
$\lambda\cdot\alpha_1\in 1+J$, and assumption (b) implies
that $\lambda=\gamma^p$ for some $\gamma\in E$. Recall that
$\underline\alpha'':=\tau_E(\gamma)\cdot\underline\alpha'=
(\gamma^{p^n}\cdot\alpha'_n~|~n\in\N)$ (proposition
\ref{prop_Teich-series}(i)); we may then replace
$\underline\alpha'$ by $\underline\alpha''$, after
which we may assume that $\alpha'_1\equiv 1\pmod{J}$.
Say that $\alpha'_0=\sum_{i=1}^kc_i\beta_i$ for some
$c_1,\dots,c_k\in E$, and set
$$
\underline\delta:=\underline\alpha'-
\sum_{i=1}^k\tau_E(c_i)\cdot\tau_E(\beta_i).
$$
Notice that $\delta_0=0$, and the images of
$\underline\alpha'$ and $\underline\delta$ agree in
$W(E/J)$. Invoking (b) again, we may find
$\gamma\in E$ such that $\gamma^p=\delta_1$, and
after replacing $\underline\alpha'$ by
$\tau_A(\gamma^{p-1})\cdot\underline\alpha'$, and
$c_i$ by $\gamma^{p-1}c_i$ for $i=1,\dots,k$, we may
further assume that $\delta_1\equiv 1\pmod{\Phi_E(J)}$.
In this case, say that $\delta_1=1+\sum^k_{i=1}d_i^p\beta_i^p$
for some $d_1,\dots,d_k\in E$; taking into account
remark \ref{rem_homogeneous-laws}(ii) and
\eqref{eq_F-V-commute} we see that
$$
\underline\delta\equiv p+p\cdot
\sum_{i=1}^k\tau_E(d_i)\cdot\tau_E(\beta_i)\pmod{p^2W(E)}.
$$
Summing up, we conclude that there exists
$\underline x\in W(E)$ such that
$p+p^2\underline x-\underline\alpha'$ lies in the ideal
generated by $(\tau_E(\beta_1),\dots,\tau_E(\beta_k))$.
Hence, $p\cdot(1+p\cdot u(\underline x))$ lies in the
ideal generated by
$(u\circ\tau_E(\beta_1),\dots,u\circ\tau_E(\beta_k))$.
Lastly, we have $1+p\cdot\underline x\in W(E)^\times$, since
$p$ is topologically nilpotent in $W(E)$, whence (c).

Obviously, (e)$\Rightarrow$(d), so we suppose that (d)
holds, and we show (e). Arguing as in the foregoing, we
may assume that $\alpha'_1\equiv 1\pmod{J}$. Next,
using (b) we find $\gamma\in E$ such that
$\gamma^p=\alpha'_1$, and after replacing $\alpha'$ by
$\tau_E(\gamma^{p-1})\cdot\underline\alpha'$, we may
even assume that $\alpha'_1\equiv 1\pmod{\Phi_E(J)}$.

Likewise, we may find $\delta\in E$ such that
$\delta^p=\alpha_1^{-1}$, and after replacing
$\underline\alpha$ by
$\tau_E(\delta)\cdot\underline\alpha$ we may assume
that $\alpha_1=1$. Now, write $\alpha'_1=1-x^p$ for
some $x\in J$, and set
$$
\underline\alpha'':=\underline\alpha'+
\tau_E(x)\cdot\underline\alpha.
$$
In light of remark \ref{rem_homogeneous-laws}(ii), we
easily see that $\alpha''_0\in J$ and
$\alpha''_1-1\in\alpha_0E$. Since $\alpha_0$ lies in
the Jacobson radical of $E$, it follows that
$\alpha''_1\in E^\times$, whence (e).

Lastly, suppose that $\cI=\underline\alpha W(E)$. In this
case, we already know that (f)$\Rightarrow$(e). Conversely,
if (e) holds, remark \ref{rem_homogeneous-laws}(ii) shows
that there exists
$\underline\beta:=(\beta_n~|~n\in\N)\in W(A)$ with
$$
\alpha_0\beta_0\in J
\qquad\text{and}\qquad
\alpha_0^p\beta_1+\alpha_1\beta_0^p\in E^\times.
$$
Then, since $\alpha_0$ is in the Jacobson radical of $E$,
we deduce that $\alpha_1\beta_0^p\in E^\times$, so also
$\beta_0\in E^\times$, and (f) follows.
\end{proof}

\sset\subsubsection{}\label{subsec_transvesal}
Additionally, often we will be dealing with pairs
$(\underline a,\cK)$ consisting of an ideal $\cK$
in a given $\F_p$-algebra $E$, and an element
$\underline a\in W(E)$ (usually $\underline a$
shall be a distinguished element); then we shall say
that $(\underline a,\cK)$ is a {\em transversal pair}
if we have
$$
W(\cK)\cap\underline a W(E)=\underline aW(\cK).
$$
This property plays a crucial role in several questions
about perfectoid rings, and the following proposition
collects a few example of such pairs. One more case
is given by theorem \ref{th_taut-two}(iv).

\begin{proposition}\label{prop_square-powers}
Let $E$ be any perfect $\F_p$-algebra, and
$\underline a:=(a_n~|~n\in\N)\in W(E)$ any element.
The following holds :
\begin{enumerate}
\item
Suppose that $a_1\in E^\times$, and let $J_1,J_2\subset E$
be any two ideals such that
\begin{enumerate}
\item
$J_1^{\La 1\Ra}E=J_1$ and $J_2^{\La 1\Ra}E=J_2$
\item
there exist $n,m\in\N$ such that $a^n_0J_1\subset J_2$
and $a_0^mJ_2\subset J_1$.
\end{enumerate}
Then the pair $(\underline a,J_1)$ is transversal if and
only if the same holds for the pair $(\underline a,J_2)$.
\item
Suppose that $E=\sum_{n\in\N}a_nE$. Then the pair
$(\underline a,I^{\lfloor r\rfloor}E)$ is transversal,
for every ideal $I\subset E$ and every $r\in\R_+$
(notation of remark {\em\ref{rem_Witt-are-f-adic}(i)}).
\item
Let $J\subset E$ be any ideal such that $J^{\La 1\Ra}E=J$.
Then the pair $(p,J)$ is transversal.
\end{enumerate}
\end{proposition}
\begin{proof}(iii) is immediate from \eqref{subsec_perfect-case}
and remark \ref{rem_Witt-are-f-adic}(iv).

(i): First, we show that the assertion holds
in the special case where $n=0$ and $m=1$. To begin
with, from lemma \ref{lem_Witt-truncate}(i) and our
assumption on $a_1$ we obtain
$$
\underline a=\tau_E(a_0)+p\underline u
\qquad
\text{for some $\underline u\in W(E)^\times$}.
$$
Notice as well that $\underline a$ is a regular element of
$W(E)$, by virtue of proposition \ref{prop_reduced-Witt}(i).
Now, suppose first that the pair $(\underline a,J_1)$
is transversal, and let $x\in W(E)$ be any element such
that $\underline a\cdot x\in W(J_2)$; then
$\tau_E(a_0)\cdot\underline a\cdot x\in W(J_1)$
(proposition \ref{prop_Teich-series}(i)), so
$\tau_A(a_0)\cdot x\in W(J_1)$, by assumption. Hence
$p\underline u\cdot x=(\underline a-\tau_A)\cdot x\in W(J_2)$,
so $px\in W(J_2)\cap pW(E)=pW(J_2)$, where the last
equality follows from (i.a) and (iii); therefore
$x\in W(J_2)$, as $p$ is regular in $W(E)$.

Likewise, suppose that $(\underline a,J_2)$ is a
transversal pair, and let $x\in W(E)$ be any element
such that $\underline a\cdot x\in W(J_1)$; then
$\underline a\cdot x\in W(J_2)$, so $x\in W(J_2)$,
and $\tau_E(a_0)\cdot x\in W(J_2)$ (again, by
proposition \ref{prop_Teich-series}(i)). It
follows that $p\underline u\cdot x=
(\underline a-\tau_E(a_0))\cdot x\in W(J_2)$, whence
$px\in W(J_2)\cap pW(E)=pW(J_2)$, so $x\in W(J_2)$, as
required.

Next, let $(J_1,J_2)$ be any pair of ideals fulfilling
conditions (i.a) and (i.b). Notice that we have
$(a_0^iJ_1)^{\La 1\Ra}E=a^i_0J_1$ for every $i\in\N$
(lemma \ref{lem_mon-fract-powers}(iii)), and the
foregoing case implies that assertion (i) holds
for the pair $(a_0^iJ_1,a_0^{i+1}J_1)$, for every $i\in\N$.
It follows that assertion (i) also holds for the
pair $(J_1,a_0^nJ_1)$; consequently we are reduced
to showing assertion (i) for the pair $(J_2,a_0^nJ_1)$,
and notice that $a^{m+n}_0J_2\subset a_0^nJ_1$.
Thus, we may replace $J_1$ by $a_0^nJ_1$, $m$ by $n+m$,
and assume from start that $J_1\subset J_2$ and
$a_0^mJ_2\subset J_1$. Now, we have
$$
a_0\cdot(J_1+a_0^iJ_2)^{\La 1\Ra}E=
(a_0J_1+a_0^{i+1}J_2)^{\La 1\Ra}E\subset
(J_1+a_0^{i+1}J_2)^{\La 1\Ra}E
\qquad
\text{for every $i\in\N$}
$$
so that -- again by the foregoing case -- we know that
assertion (i) holds for the pair
$$
((J_1+a_0^iJ_2)^{\La 1\Ra}E,(J'+a_0^{i+1}J_2)^{\La 1\Ra}E)
\qquad
\text{for every $i\in\N$}
$$
and then it also holds for the pair
$((J_1+a_0^mJ_2)^{\La 1\Ra}E,(J_1+J_2)^{\La 1\Ra}E)$; but
$(J_1+a_0^mJ_2)^{\La 1\Ra}E=J_1$ and
$(J_1+J_2)^{\La 1\Ra}E=J_2$, whence the contention.

(ii): Indeed, let $\underline b\in W(A)$ be any
element, and suppose that
$\underline a\cdot\underline b\in W(I^{\lfloor r\rfloor}E)$. 
Define $\cR_E$ as in lemma \ref{lem_new-powers}; by lemma
\ref{lem_new-powers}(iv) and proposition
\ref{prop_semi-norm-on-W}(iii), we deduce that
$$
|\underline a|_1\cdot|\underline b|_1=
|\underline a\cdot\underline b|_1\leq|I|^r
\qquad
\text{for every $|\cdot|$ in $\cR_E$}.
$$
But since the system $\{a_n~|~n\in\N\}$ generates $E$,
it is easily seen that $|\underline a|_1=1$ for every
$|\cdot|$ in $\cR_E$. We conclude that
$|\underline b|_1\leq|I|^r$ for every such $|\cdot|$,
and then the assertion follows, again by invoking
lemma \ref{lem_new-powers}(iv).
\end{proof}

\begin{corollary}\label{cor_square-powers}
Let $E$ be a perfect $\F_p$-algebra, $b_\bullet:=(b_1,\dots,b_n)$
a finite system of elements of $E$, and $\cI\subset W(E)$ the
ideal generated by $(\tau_E(b_1)\dots,\tau_E(b_n))$. Let also
$\underline a:=(a_n~|~n\in\N)\in W(E)$ be any element, and
suppose that :
\begin{enumerate}
\alphaenu
\item
$E=\sum_{n\in\N}a_nE$.
\item
$W(E)$ is complete and separated for its $\cI$-adic
topology.
\end{enumerate}
Then $\underline aW(E)$ is a closed ideal for the
$\cI$-adic topology of\/ $W(E)$.
\end{corollary}
\begin{proof} Let $I\subset E$ be the ideal generated
by $b_\bullet$. Taking into account proposition
\ref{prop_morel}(i) and lemma
\ref{lem_mon-fract-powers}(iv), we easily see that
the $\cI$-adic topology on $W(E)$ agrees with the
linear topology defined by the cofiltered system of
ideals $(W(I^{\lfloor q\rfloor})~|~q\in\Q_+)$. Then, propositions
\ref{prop_square-powers}(ii) \ref{prop_reduced-Witt}(i)
imply that scalar multiplication by $\underline a$ is a
$W(E)$-linear isomorphism
$$
W(E)\isom\underline aW(E)
$$
that identifies the $\cI$-adic topology of
$\underline aW(E)$ with the one induced by the
$\cI$-adic topology of $W(E)$. Since $W(E)$
is complete and separated for its $\cI$-adic topology,
it follows that $\underline aW(E)$ is complete
and separated for the subspace topology induced from
the $\cI$-adic topology of $W(E)$, and the
contention follows.
\end{proof}

\subsection{P-rings}
In this section we make a preliminary study of an auxiliary
class of topological rings containing the perfectoid rings
that shall be introduced in section \ref{sec_now-perfectoid}.

\begin{definition}\label{def_perfectoid}
A complete and separated topological ring $(A,\cT)$ is
called a {\em P-ring} if there exist a prime integer
$p$ and an ideal $I\subset A$ such that the following
conditions hold :
\begin{enumerate}
\alphaenu
\item
$I$ is finitely generated and $\cT$ agrees with the
$I$-adic topology.
\item
$pA\subset I^2$ and the Frobenius endomorphism of $A/I^2$
is a surjection 
$$
\Phi_{A/I^2}:A/I^2\to A/I^2
\qquad
a\mapsto a^p.
$$
\end{enumerate}
Any ideal $I$ of $A$ fulfilling conditions (a) and (b)
is called an {\em ideal of definition} of $A$. Moreover,
we say that $I$ is {\em special}, if the following holds :
\begin{enumerate}
\alphaenu
\addenu\addenu
\item
There exists a finite set of generators $a_1,\dots,a_n$
of $I$ such that $pA\subset I^{(p)}:=\sum_{i=1}^nAa^p_i$.
\end{enumerate}
\end{definition}

It turns out that if $I$ is a special ideal of definition
of the P-ring $A$, the ideal $I^{(p)}$ depends only on $I$
(and not on the choice of a system of generators for $I$),
so this notation is not abusive. Indeed, we have more
generally :

\begin{lemma}\label{lem_special-ideals}
Let $A$ be any ring, $p$ a prime integer, and
$(a_1,\dots,a_n),(b_1,\dots,b_m)$ two sequences
of elements of $A$. Suppose that
$$
I:=\sum_{i=1}^nAa_i\subset J:=\sum_{i=1}^mAb_i
\qquad\text{and}\qquad
p\in I^{(p)}:=\sum_{i=1}^nAa_i^p.
$$
Then :
\begin{enumerate}
\item
$I^{(p)}\subset J^{(p)}:=\sum_{i=1}^mAb^p_i$.
\item
Especially, if $I=J$, we have $I^{(p)}=J^{(p)}$.
\end{enumerate}
\end{lemma}
\begin{proof} It is easily seen that (i)$\Rightarrow$(ii),
so we need only check (i). To this aim, it suffices to
show that $I^{(p)}_\fm\subset J^{(p)}_\fm$ for every maximal
ideal $\fm\subset A$, so we may replace $A$ by $A_\fm$ and
assume from start that $A$ is local. Now, suppose first
that $J=A$; then it is easily seen that $J^{(p)}=A$ as
well, so the assertion is clear in this case. Hence, we
may also assume that $J$ is contained in the maximal ideal
of $A$. By assumption, for every $i=1,\dots,n$ we may write
$a_i=\sum_{j=1}^mx_{ij}b_j$ for certain $x_{i1},\dots,x_{im}\in A$;
it follows easily that $a_i^p-\sum_{j=1}^mx_{ij}^pb_j^p\in pJ^p$
(details left to the reader). Summing up, we see that
$I^{(p)}\subset J^{(p)}+pJ^p\subset J^{(p)}+I^{(p)}J^{(p)}$,
and the assertion follows by Nakayama's lemma.
\end{proof}

\begin{lemma}\label{lem_perfectoid}
Let $A$ be any P-ring, and $I$ any ideal of
definition of $A$. We have :
\begin{enumerate}
\item
If $A\neq 0$, there exists a unique prime integer $p$,
independent of $I$, such that condition {\em (b)} of
definition {\em\ref{def_perfectoid}} is fulfilled.
\item
$I$ is contained in the Jacobson radical of $A$.
\item
There exists a special ideal of definition $J$ such
that $J^{(p)}=I$.
\item
For $p$ as in {\em (i)}, the Frobenius endomorphism
$\Phi_{A/pA}:A/pA\to A/pA$ is a surjective ring homomorphism,
and there exist $\pi\in A$ and $u\in A^\times$ such that
$p=u\pi^p$.
\item
For $p$ as in {\em (i)}, the ring $A$ is $p$-adically
complete and separated.
\end{enumerate}
\end{lemma}
\begin{proof}(ii) is standard : if $a\in I$, then
$1-ab$ is invertible in $A$ for every $b\in A$, since
$A$ is $I$-adically complete and separated (the inverse
is given by the convergent series $\sum_{k\in\N}(ab)^n$);
the claim follows.

(i): Let $J$ be another ideal of definition of $A$,
and $p,q\in\Z$ two prime integers with $p\in I^2$
and $q\in J^2$. By (ii), both $p$ and $q$ lie in the
Jacobson radical of $A$; as $A\neq 0$, we get $p=q$.

(iii): Let $(f_\lambda~|~\lambda\in\Lambda)$ be a
finite system of generators (indexed by the finite
set $\Lambda$) for $I$; by assumption, we may find
$g_\lambda\in A$ such that $f_\lambda-g_\lambda^p\in I^2$
for every $\lambda\in\Lambda$. From (ii) and Nakayama's
lemma, we deduce that the system
$(g^p_\lambda~|~\lambda\in\Lambda)$ generates $I$, so
we can take for $J$ the ideal generated by the system
$(g_\lambda~|~\lambda\in\Lambda)$.

(iv): Let $(g_\lambda~|~\lambda\in\Lambda)$ be as in
the proof of (iii). For any $a\in A$ we construct
inductively a sequence $(b_n~|~n\in\N)$ of elements
of $A$, such that
$$
a-\sum_{k=0}^n b^p_k\in I^{n+1}
\qquad\text{and}\qquad
b_n\in I^n
\qquad
\text{for every $n\in\N$}.
$$
Indeed, for $n=0$, since $\Phi_{A/I}$ is surjective
we may find $b_0\in A$ such that $a-b_0^p\in I$. Next,
suppose that $n\geq 0$ and $b_0,\dots,b_n\in A$ are
already given such that
$c_n:=a-\sum_{k=0}^n b^p_k\in I^{n+1}$; we may find a system
$(d_{\underline\lambda}~|~\underline\lambda\in\Lambda^{n+1})$
of elements of $A$ such that
$$
c_n=\sum_{\underline\lambda\in\Lambda^{n+1}}
g^p_{\underline\lambda}d_{\underline\lambda}
\qquad\text{where}\qquad
g_{\underline\lambda}:=\prod^{n+1}_{i=1}g_{\lambda_i}
\quad
\text{for every
$\underline\lambda:=(\lambda_1,\dots,\lambda_{n+1})\in\Lambda^{n+1}$}.
$$
In turn, for each $\underline\lambda\in\Lambda^{n+1}$ we may write
$d_{\underline\lambda}=e^p_{\underline\lambda}+e'_{\underline\lambda}$
for elements $e_{\underline\lambda}\in A$ and
$e'_{\underline\lambda}\in I$, so that
$$
c_n=c'_n+c''_n
\qquad\text{where}\quad
c'_n:=\sum_{\underline\lambda\in\Lambda^{n+1}}
g^p_{\underline\lambda}e^p_{\underline\lambda}
\qquad
c''_n:=\sum_{\underline\lambda\in\Lambda^{n+1}}
g^p_{\underline\lambda}e'_{\underline\lambda}\in I^{n+2}.
$$
Thus, $b_{n+1}:=\sum_{\underline\lambda\in\Lambda^{n+1}}
g_{\underline\lambda}e_{\underline\lambda}\in I^{n+1}$,
and since $p\in I^2$, a simple computation shows that
$c'_n-b_{n+1}^p\in I^{n+2}$, so finally
$a-\sum^{n+1}_{k=0}b_k^p\in I^{n+2}$. Notice that the
series $\sum_{k\in\N}b^p_k$ converges to $a$ in the
$I$-adic topology of $A$, and set $b:=\sum_{k\in\N}b_k$;
it follows easily that $a-b^p\in pA$, whence the surjectivity
of $\Phi_{A/pA}$. Lastly, by assumption we may write
$p=\sum_{i=1}^ra_ib_i$ for finitely many elements
$a_1,b_1,\dots,a_r,b_r\in I$; by the foregoing, for every
$i=1,\dots,r$ we may further find $f_i,g_i,f'_i,g'_i\in A$
such that $a_i=f_i^p+pf'_i$, $b_i=g_i^p+pg'_i$. Then
$f_i^p,g_i^p\in I$, and (ii) implies that $f_i$ and $g_i$
lie in the Jacobson radical $\rad(A)$ of $A$, for every
$i\leq r$. Set $\pi:=\sum_{i=1}^rf_ig_i$; summing up, we
find $p-\pi^p=pc$ for some $c\in\rad(A)$, so
$p\cdot(1-c)=\pi^p$, and to conclude it suffices to
remark that $1-c\in A^\times$.

(v) is a special case of lemma \ref{lem_fontaine}.
\end{proof}

\begin{remark}\label{rem_why-only-now}
Suppose that $A$ is a complete and separated topological
$\F_p$-algebra. Then it follows easily from lemma
\ref{lem_perfectoid}(iv) that $A$ is a P-ring if and only if
its Frobenius endomorphism $\Phi_A$ is surjective, and
there exists a finitely generated ideal $I\subset A$ such
that the topology of $A$ agrees with the $I$-adic topology.
\end{remark}

\begin{proposition}\label{prop_P-rings}
Let $A$ be a ring, $p$ a prime integer, $I\subset A$ a
finitely generated ideal such that $pA\subset I$, and
denote by $\cT_I$ (resp. $\cT_p$) the $I$-adic (resp.
$p$-adic) topology on $A$. We have :
\begin{enumerate}
\item
If $(A,\cT_I)$ is a P-ring, the same holds for $(A,\cT_p)$.
\item
If $(A,\cT_p)$ is a P-ring, the same holds for the
completion $(A^\wedge_I,\cT_I^\wedge)$ of $(A,\cT_I)$.
\end{enumerate}
\end{proposition}
\begin{proof} (i) follows immediately from lemma
\ref{lem_perfectoid}(iv,v).

(ii): Let $(f_1,\dots,f_k)$ be a finite system of
generators of $I$; by lemma \ref{lem_perfectoid}(iv)
there exists $\pi\in A$ such that $pA=\pi^pA$, and
the Frobenius endomorphism $\Phi_{A/pA}$ of $A/pA$ is
surjective, so for every $i=1,\dots,k$ we may find
$g_i\in A$ such that $f_i-g_i^p\in pA$, and we let
$J\subset A$ be the ideal generated by the finite
system $(g_1,\dots,g_k,\pi)$. Set $N:=(p-1)(k+1)+1$,
and notice that $J^N\subset I\subset J^p$. Especially,
$\cT_I$ agrees with the $J$-adic topology on $A$, and
moreover $p\in J^2$. Furthermore, $\cT_I^\wedge$ agrees
with the $JA^\wedge_I$-adic topology on $A^\wedge_I$
(remark \ref{rem_completion-of-topring}(iv) and
lemma \ref{lem_still-c-adic}(iv)), and the natural map
$A/J^2\to A^\wedge_I/J^2A^\wedge_I$ is an isomorphism
(remark \ref{rem_completion-of-topring}(ii,iv));
since $\Phi_{A/pA}$ is surjective by assumption, the
same holds for $\Phi_{A/J^2}$, and the assertion follows.
\end{proof}

\sset\subsubsection{}\label{subsec_back-to-E-A}
{\em Henceforth we fix a prime integer $p$, and we shall
assume that the topology of every P-ring that appears
throughout the rest of this chapter, is coarser than the
$p$-adic topology}.
Let $(A,\cT)$ be any P-ring and $I\subset A$ an ideal of
definition. Notice that if $I$ is special (definition
\ref{def_perfectoid}), the Frobenius endomorphism
$\Phi_{A/I}$ of $A/I$ factors uniquely as the composition
of a morphism of topological rings
$$
\bar\Phi_{A/I}:A/I\to A/I^{(p)}
$$
and the natural projection $A/I^{(p)}\to A/I$.

\begin{lemma}\label{lem_was-third-cond}
In the situation of \eqref{subsec_back-to-E-A}, we have :
\begin{enumerate}
\item
The maps $u_A:\bA(A)\to A$ and $\bar u_{A/pA}:\bE(A)\to A/pA$
are open and surjective.
\item
$\Ker\,u_A$ contains a distinguished element of\/ $\bA(A)$
(see definition {\em\ref{def_distinguished}}).
\item
If $\underline\alpha:=(\alpha_n~|~n\in\N)$ is any
distinguished element in $\Ker\,u_A$, we have
$$
\bar u_A(\alpha_0)=pu
\qquad
\text{for some $u\in A^\times$}.
$$
\end{enumerate}
\end{lemma}
\begin{proof} (i): Let $I\subset A$ be any ideal of
definition of $A$, so that the topology of $A/pA$ agrees
with the $I/pA$-adic topology; by lemma
\ref{lem_perfectoid}(iv) we know already that $\Phi_{A/pA}$
is surjective, and then it is easily seen that there exists
$n\in\N$ such that $(I/pA)^n\subset\Phi_{A/pA}(I)$. Especially,
$\Phi_{A/pA}$ is an open map; then, the assertion for
$\bar u_{A/pA}$ follows from remark \ref{rem_topology-of-E}(iii).
By virtue of \cite[Th.8.4]{Mat}, we deduce the surjectivity of
$u_A$. Next, let  $(a_1,\dots,a_r)$ be a finite system of
generators for $I$, and denote by $\bar a_i$ the image of
$a_i$ in $A/pA$, for every $i=1,\dots,k$. Since $\bar u_{A/pA}$
is surjective, we may find
$\bar\alpha_1,\dots,\bar\alpha_k\in\bE(A)$ such that
$\bar u_{A/pA}(\bar\alpha_i)=\bar a_i$ for $i=1,\dots,k$,
and we set $\alpha_i:=\tau_A(\bar\alpha_i)$ for $i=1,\dots,k$
(notation of \eqref{sec_def-A-tilde}).
Since $I$ is topologically nilpotent in $A$, it is
easily seen that $\bar\alpha_i$ is topologically
nilpotent in $\bE(A)$, and since $\tau_A$ is continuous,
we deduce that $\alpha_i$ is topologically nilpotent
in $\bA(A)$, for every $i=1,\dots,k$. Let $\cI\subset\bA(A)$
be the ideal generated by the system
$(\alpha_1,\dots,\alpha_k)$; it follows easily that
$\cI$ is topologically nilpotent as well. On the
other hand, by construction we have
$a_i-u_A(\alpha_i)\in pA\subset I^2$ for every
$i=1,\dots,k$, so that $u_A(\cI)=I$, by Nakayama's lemma.
Lastly, let $\cJ\subset\bA(A)$ be any open ideal; then
$\cI^n\subset\cJ$ for every sufficiently large $n\in\N$,
so $I^n\subset u_A(\cI)$, and the assertion follows.

(ii): We notice the following :

\begin{claim}\label{cl_was-assertion-ii}
There exist $w\in A^\times$ and $t\in\bE(A)$ such
that $p=w\cdot\bar u_A(t)$ in $A$.
\end{claim}
\begin{pfclaim} By lemma \ref{lem_perfectoid}(iv), we may
write $p=v\cdot x^p$ for some $v\in A^\times$ and $x\in A$.
Since $\bar u_{A/pA}$ is surjective, we may find
$s\in\bE(A)$ such that $x\equiv\bar u_A(s)\pmod{pA}$.
In view of lemma \ref{lem_basic-cong}(i) we deduce that
$p\equiv v\cdot\bar u_A(s)^p\pmod{p^2A}$, {\em i.e.}
$p=v\cdot\bar u_A(t)+p^2b$ for $t:=s^p$ and some $b\in A$;
thus $p\cdot(1-pb)=v\cdot\bar u_A(t)$.
But since $A$ is $p$-adically complete and separated
(lemma \ref{lem_perfectoid}(v)), $1-pb\in A^\times$
(remark \ref{rem_someth-on-bdd-in-Z-lin}(v)), so we
get the sought identity with $w:=v\cdot(1-pb)^{-1}$.
\end{pfclaim}

Now, take $w$ and $t$ as in claim \ref{cl_was-assertion-ii}.
Since $u_A$ is surjective, we may find
$\beta:=(b_n~|~n\in\N)\in\bA(A)$ with $u_A(\beta)=w$.
By proposition \ref{prop_Teich-series}(i) we have
$\beta\cdot\tau_A(t)=(t^{p^n}b_n~|~n\in\N)$, and
\eqref{eq_this-is-p} and remark \ref{rem_homogeneous-laws}(ii)
show that
$$
\gamma:=\beta\cdot\tau_A(t)-p=(tb_0,t^pb_1,\dots)-(0,1,\dots)
=(tb_0,t^pb_1-1,\dots).
$$
But notice that $t$ is topologically nilpotent in $\bE(A)$;
since the topology of $\bE(A)$ is linear, it follows that
$tb_0$ and $t^pb_1$ are topologically nilpotent
in $\bE(A)$, and since $\bE(A)$ is complete and separated,
$1-t^pb_1\in\bE(A)^\times$, so $\gamma$ is distinguished,
and it lies in the kernel of $u_A$, due to lemma
\ref{lem_drop-conditions}(iv).

(iii): Indeed, from \eqref{eq_F-V-commute} and lemma
\ref{lem_Witt-truncate}(i) we get
$\underline\alpha=\tau_A(\alpha_0)+p\cdot\underline\beta$,
where $\beta:=(\alpha_{n+1}^{1/p}~|~n\in\N)$, and since
$\alpha_1\in\bE(A)^\times$, we have $\underline\beta\in\bA(A)^\times$
(proposition \ref{prop_Witt-is-complete}(iii)); now,
$\bar u_A(\alpha_0)=-p\cdot u_A(\underline\beta)$, whence
the assertion, in light of lemma \ref{lem_drop-conditions}(iv).
\end{proof}

\begin{proposition}\label{prop_third-equiv-cond}
In the situation of \eqref{subsec_back-to-E-A}, 
the following conditions are equivalent :
\begin{enumerate}
\alphaenu
\item
$\bE(A)$ is a P-ring.
\item
There exists an ideal of definition $I$ of $A$
such that $\Ker\,\Phi_{A/I}$ is a finitely generated
ideal of $A/I$.
\item
There exists a special ideal of definition $J$ of
$A$ such that\/ $\bar\Phi_{A/J}$ is an isomorphism.
\end{enumerate}
\end{proposition}
\begin{proof} For any ideal of definition $I$ of $A$,
recall that $\bE(A/I)$ is the limit of the inverse system
$(A_n~|~n\in\N)$ with $A_n:=A/I$ for every $n\in\N$, and
with transition maps $\phi_n:A_{n+1}\to A_n$ given by the
Frobenius endomorphism $\Phi_{A/I}$ of $A/I$. For every
$n\in\N$ we shall denote by
$$
\bar u_n:\bE(A)\to A_n
$$
the composition of the canonical identification
$\bE(A)\isom\bE(A/I)$ as in \eqref{sec_def-A-tilde},
with the natural projection $\bE(A/I)\to A_n$.
Especially $\bar u_0=\bar u_{A/I}$; also, each $\bar u_n$
is surjective, since the same holds for $\Phi_{A/I}$.
Set $\cJ_0:=\Ker\,\bar u_{A/I}$, and let $\Phi_{\bE(A)}$
be the Frobenius endomorphism of $\bE(A)$; then the
family of ideals
$$
\cJ_n:=\Phi^n_{\bE(A)}(\cJ_0)=\Ker\,\bar u_n
\qquad
\text{for every $n\in\N$}
$$
is a fundamental system of open neighborhoods of $0$ for
the topology $\cT_{\bE(A)}$ of $\bE(A)$. Especially,
every element of $\cJ_0$ is topologically nilpotent, and
since $\bE(A)$ is complete and separated (remark
\ref{rem_topology-of-E}(i)), we easily deduce that
$\cJ_0$ is contained in the Jacobson radical of $\bE(A)$
(cp. the proof of lemma \ref{lem_perfectoid}(ii)).

(c)$\Rightarrow$(a): Let $J$ be as in (c), and take $I:=J$
in the foregoing. Pick a finite system $(a_1,\dots,a_r)$ of
generators of $J$, and let $\pi:A\to A/J^{(p)}$ be the projection.
A simple inspection of the definition yields the identity :
$$
\bar u_{A/J^{(p)}}=\bar\Phi_{A/J}\circ\bar u_1.
$$
Since $\bar\Phi_{A/J}$ is an isomorphism, we deduce that
$\bar u_{A/J^{(p)}}$ is a surjection and we may find
$$
\alpha_1,\dots,\alpha_r\in\bE(A)
\quad\text{with}\quad
\bar u_{A/J^{(p)}}(\alpha_i)=\pi(a_i)
\qquad
\text{for every $i=1,\dots,r$}.
$$
Let $\cJ'_0\subset\bE(A)$ be the ideal generated by
the system $(\alpha_1,\dots,\alpha_r)$, and set
$$
\cJ'_n:=\Phi_{\bE(A)}^n(\cJ'_0)
\qquad\text{for every $n\in\N$}.
$$

\begin{claim}\label{cl_drop-second}
$\cJ'_0=\cJ_0$, and the $\cJ_0$-adic topology on $\bE(A)$
agrees with $\cT_{\bE(A)}$.
\end{claim}
\begin{pfclaim} A simple inspection shows that
$\cJ'_0\subset\cJ_0$, and consequently $\cJ'_n\subset\cJ_n$
for every $n\in\N$. Moreover, $\bar u_{A/J^{(p)}}(\cJ'_0)=J/J^{(p)}$,
and therefore
$$
\bar u_1(\cJ'_0)=\bar\Phi_{A/J}^{\ -1}\circ\bar u_{A/J^{(p)}}(\cJ'_0)
=\bar\Phi_{A/J}^{\ -1}(J/J^{(p)})=\Ker\,\phi_0=\cJ_0/\cJ_1.
$$
Especially, the natural map
$\cJ'_0/\cJ'_1\to\cJ_0/\cJ_1$ is surjective, and since
$\Phi_{\bE(A)}$ is an automorphism, we deduce that the
induced map $\cJ'_n/\cJ'_{n+1}\to\cJ_n/\cJ_{n+1}$ is
surjective for every $n\in\N$. The filtration
$(\cJ_n~|~n\in\N)$ is separated on $\bE(A)$, so the
identity $\cJ_0=\cJ'_0$ will follow from
\cite[Ch.III, \S2, n.8, Cor.2]{BouAC}, after we show that
$\bE(A)$ is complete for the topology $\cT$ defined by
the filtration $(\cJ'_n~|~n\in\N)$. However, since
$\cJ'_0$ is finitely generated, for every $n\in\N$
there exists $k(n)\in\N$ such that
$\cJ'_{k(n)}\subset\cJ_0^{\prime k(n)}\subset\cJ'_n$,
so the topology $\cT$ agrees with the $\cJ'_0$-adic
topology, and $\cJ'_0$ is topologically nilpotent
for the topology $\cT_{\bE(A)}$, since it is
contained in $\cJ_0$. Then, the assertion follows
from lemma \ref{lem_fontaine}. By the same token,
we see that the $\cJ_0$-adic topology agrees with
$\cT_{\bE(A)}$.
\end{pfclaim}

The assertion follows immediately from claim
\ref{cl_drop-second} and remark \ref{rem_why-only-now}.

(a)$\Rightarrow$(b): Let $I\subset A$ and
$\cI\subset\bE(A)$ be two ideals of definition;
then there exists $n\in\N$ such that $\cJ_n\subset\cI$.
Clearly $\cI':=\Phi_{\bE(A)}^{-n}(\cI)$ is still an
ideal of definition, and $\cJ_0\subset\cI'$.
After replacing $\cI$ by $\cI'$, we may therefore assume
that $\cJ_0\subset\cI$. Write $\bar u_{A/I}(\cI)=I'/I$ for
an ideal $I'\subset A$; we have :

\begin{claim}\label{cl_special-claim}
$\cI=\Ker\,\bar u_{A/I'}$ and $I'$ is an ideal of definition
of $A$.
\end{claim}
\begin{pfclaim} We have just remarked that
$\Ker\,\bar u_{A/I}\subset\cI$, and the stated identity
is an easy consequence. Next, since both $\cI$ and $I$ are
finitely generated, the same holds for $I'$, and obviously
$pA\subset I'^2$. Also $\cI$ is topologically nilpotent
in $\bE(A)$, and the topology of $A/I$ is discrete, hence
$I'/I$ is a nilpotent ideal of $A/I$, and it follows
easily that the $I$-adic topology on $A$ agrees with
the $I'$-adic topology.
\end{pfclaim}

Due to claim \ref{cl_special-claim}, we may replace
$I$ by $I'$, and assume from start that $\cI=\cJ_0$.
We may write $\bar u_1(\cJ_0)=J_1/I$ for an ideal
$J_1\subset A$ containing $I$, and since $\cJ_0$
and $I$ are finitely generated, the same holds for
$J_1$. It suffices now to notice that
\set\begin{equation}\label{eq_ker-of-phi_0}
J_1/I=\Ker\,\phi_0.
\end{equation}

(b)$\Rightarrow$(c): By assumption, there exists
a finitely generated ideal $J_1\subset A$ containing
$I$ and such that \eqref{eq_ker-of-phi_0} holds.
Furthermore, by lemma \ref{lem_perfectoid}(iii) we may
write $I=J'{}^{(p)}$ for some ideal of definition
$J'\subset A$. Set $J:=J'+J_1$. From \eqref{eq_ker-of-phi_0}
we get $J_1^{(p)}\subset I$, therefore
$J^{(p)}=J'{}^{(p)}+J_1^{(p)}=I$. We conclude that
$\Phi_{A/J}$ induces an isomorphism
$A/J_1\isom A/I=A/J^{(p)}$; but the latter factors
through $\bar\Phi_{A/J}$ and the projection
$A/J_1\to A/J$, so $\bar\Phi_{A/J}$ is an isomorphism
as well (and $J=J_1$). Lastly, since $\cJ_0$ is
topologically nilpotent, $J_1/I$ is a nilpotent
ideal, so the same holds for $J/J'$, and
consequently the $J$-adic topology agrees with the
$J'$-adic topology; so $J$ is an ideal of definition,
whence (b).
\end{proof}

Let $A$ be any P-ring; by lemma \ref{lem_was-third-cond}(ii),
the kernel of $u_A$ contains a distinguished element
$\underline a$, and we shall see later that the case
where actually $\underline a$ generates the ideal
$\Ker\,u_A$ is the most interesting for our purposes.
For the moment, let us just notice that, if the latter
condition holds, the datum $(E:=\bE(A),\underline a)$
suffices to recover $A$ : indeed, lemma
\ref{lem_was-third-cond}(i) says that in that case
$u_A$ induces an isomorphism of topological rings
$W(E)/\underline aW(E)\isom A$ (where $W(E)/\underline aW(E)$
is endowed with the quotient topology induced from
$W(E)$ via the natural projection). It is then convenient
to insert hereafter a preliminary study of pairs
$(E,\underline a)$ of this type.

\sset\subsubsection{}\label{subsec_several-tops}
Let $(E,\cT_E)$ be a perfect topological $\F_p$-algebra,
$\underline a:=(a_n~|~n\in\N)\in W(E)$ any element, and
suppose that $\cT_E$ agrees with the $a_0$-adic topology.
Denote by $E^\wedge$ the separated completion of $E$, and
let  also $E_d$ (resp. $E^\wedge_d$) be the ring $E$ (resp.
$E^\wedge$) endowed with its discrete topology. To ease
notation, set
$$
A(E):=W(E)/\underline aW(E)
$$
which we view as a topological ring, with the quotient
topology induced from $W(E)$. Define likewise the topological
rings $A(E_d)$, $A(E^\wedge)$ and $A(E^\wedge_d)$.
The natural diagram of continuous ring homomorphisms
$$
\xymatrix{
E_d \ar[r] \ar[d] & E^\wedge_d \ar[d] \\
E \ar[r] & E^\wedge
}$$
induces a commutative diagram of topological rings
\set\begin{equation}\label{eq_many-quotients}
{\diagram
A(E_d) \ar[rr] \ar[d] & & A(E^\wedge_d)
\ar[d] \\
A(E) \ar[rr]^-{j_A} & & A(E^\wedge).
\enddiagram}
\end{equation}
Moreover, notice that there are natural isomorphisms
of topological rings
\set\begin{equation}\label{eq_mod-out-by-p}
A(E)\otimes_\Z\F_p\isom E/a_0E
\qquad
A(E^\wedge)\otimes_\Z\F_p\isom E^\wedge/a_0E^\wedge
\end{equation}
(where, again the sources and the targets are endowed
with the quotient topologies induced from $A(E)$,
$E$, $A(E^\wedge)$, and respectively $E^\wedge$). The
composition of the first of these maps with the projection
$A(E)\to A(E)\otimes_\Z\F_p$ is the continuous ring
homomorphism
$$
\pi_A:A(E)\to E/a_0E
\qquad
((b_n~|~n\in\N)\mod{\underline aW(E)})\mapsto
(b_0\mod{a_0E})
$$
and likewise we can describe the corresponding map
$\pi^\wedge_A:A(E^\wedge)\to E^\wedge/a_0E^\wedge$. Lastly,
denote
$$
\beta:\bE(A(E))\to E^\wedge
$$
the morphism of topological monoids obtained by
composing $\bE(\pi_A)$ with the isomorphism of topological
rings $\bE(E/a_0E)\isom E^\wedge$ provided by corollary
\ref{cor_E-and-completion}.

\begin{proposition}\label{prop_special-case}
In the situation of \eqref{subsec_several-tops}, suppose
that $a_0E+a_1E=E$. Then the following holds :
\begin{enumerate}
\item
The image of $\underline a$ is a distinguished element
of\/ $W(E^\wedge)$.
\item
All the arrows of \eqref{eq_many-quotients} are
isomorphisms of topological rings.
\item
$A(E)$ is a P-ring, and its topology agrees with
the $p$-adic topology.
\item
$\beta$ is an isomorphism of topological rings (for
the topological ring structure on $A(E)$ provided by
\eqref{sec_def-A-tilde}), and we have a commutative
diagram of topological rings
$$
\xymatrix{
\bA(A(E)) \ar[rr]^-{W(\beta)}
\ar[d]_{u_{A(E)}} & &
W(E^\wedge) \ar[d]^{\pi^\wedge_W} \\
A(E) \ar[rr]^-{j_A} & & A(E^\wedge)
}$$
where $\pi^\wedge_W$ is the natural projection.
\end{enumerate}
\end{proposition}
\begin{proof}(i): By remark \ref{rem_someth-on-bdd-in-Z-lin}(v),
the image of $a_0$ lies in the Jacobson radical of $E^\wedge$, and
since $a_0E+a_1E=E$, it follows easily that $a_1\in E^\times$
(details left to the reader), whence the assertion.

(ii): Let $\cI\subset W(E)$ (resp.
$\cI^\wedge\subset W(E^\wedge)$) be the ideal
generated by the system $(p,\tau_E(a_0))$ (resp.
$(p,\tau_{E^\wedge}(a_0))$); proposition \ref{prop_morel}(ii)
says that the topology of $W(E)$ (resp. of
$W(E^\wedge)$) agrees with the $\cI$-adic
(resp. $\cI^\wedge$-adic) topology, whereas the topologies
of both $W(E_d)$ and $W(E^\wedge_d)$ are $p$-adic. Let
$\bar\cI$ (resp. $\bar\cI{}^\wedge$) be the image of
$\cI$ in $A(E)$ (resp. in $A(E^\wedge)$); it follows
easily that the topology of $A(E)$ (resp. of $A(E^\wedge)$)
agrees with the $\bar\cI$-adic (resp. $\bar\cI{}^\wedge$-adic)
topology, whereas the topologies of both $A(E_d)$ and
$A(E^\wedge_d)$ are still $p$-adic.
However, set $\underline b:=(a_1^{1/p},a_2^{1/p},\dots)$,
and notice that
\set\begin{equation}\label{eq_combine-later}
\underline a=\tau_E(a_0)+p\cdot\underline b
\qquad
\text{in $W(E)$}
\end{equation}
(\eqref{eq_F-V-commute} and lemma \ref{lem_Witt-truncate}(i)),
which implies that $\bar\cI=pA(E)$. Likewise,
$\bar\cI{}^\wedge=pA(E^\wedge)$.
Summing up, this proves already that the vertical arrows
of \eqref{eq_many-quotients} are both isomorphisms.

Next, proposition \ref{prop_reduced-Witt}(v) implies
that $A(E)$ is complete and separated for
the $p$-adic topology, hence also for the $\cI$-adic
topology. Recall as well that the family of ideals
$$
\Ker\,(W(E)\to W_n(E/a_0^kE))
\qquad
\text{for every $k,n\in\N$}
$$
is a fundamental system of open neighborhoods of
$0\in W(E)$ for the $\cI$-adic topology, therefore
the natural map
$$
A(E)\to\lim_{k,n\in\N}\,
W_n(E/a_0^kE)/\underline aW_n(E/a_0^kE)
$$
is an isomorphism of topological rings. By the same
token, the same description applies to $A(E^\wedge)$,
so we conclude that also the horizontal arrows of
\eqref{eq_many-quotients} are isomorphisms.

(iii): We know already that the topology of $A(E)$ is
complete and separated, and agrees with the $p$-adic
topology. On the other hand, notice that $\underline b$
is invertible in $W(E^\wedge)$, by (i) and proposition
\ref{prop_Witt-is-complete}(iii); since the images of
$\tau_{E^\wedge}(a_0)$ and $-p\cdot\underline b$ agree
in $A(E^\wedge)$, we deduce that the $p$-adic filtration
agrees with the $\tau_E(a_0)$-adic filtration on
$A(E^\wedge)$, and therefore also on $A(E)$,
since we already know that $j_A$ is an isomorphism.
Set $J:=\tau_E(a^{1/p}_0)\cdot A(E)$; it follows
that the topology of $A(E)$ agrees with the $J$-adic
topology, and moreover $p\in J^p\subset J^2$. Lastly,
the Frobenius endomorphism is surjective on $A(E)/J^2$,
due to the isomorphism \eqref{eq_mod-out-by-p}, whence
the assertion.

(iv): From (iii) and the discussion of
\eqref{sec_def-A-tilde} we see that $\beta$ is an
isomorphism of topological rings; also, under the identifications
\eqref{eq_mod-out-by-p}, the morphism $j_A\otimes_\Z\F_p$
corresponds to the natural map
$\iota:E/a_0E\isom E^\wedge/a_0E^\wedge$. Therefore we have
$$
\pi^\wedge_A\circ j_A\circ u_{A(E)}=
\iota\circ\pi_A\circ u_{A(E)}=
\iota\circ \pi_A\circ\bar u_{A(E)}\circ\bomega_0
$$
where $\bomega_0:\bA(A(E))\to\bE(A(E))$ is the $0$-th ghost
component map. Then, since $A(E^\wedge)$ is complete and
separated for the $p$-adic topology, proposition
\ref{prop_lift-Witt}(iii) reduces to checking

\begin{claim}
$\pi^\wedge_A\circ\pi^\wedge_W\circ W(\beta)=
\iota\circ\pi_A\circ\bar u_{A(E)}\circ\bomega_0$.
\end{claim}
\begin{pfclaim}[] Notice the identities :
$$
\beta\circ\bomega_0=\bomega^\wedge_0\circ W(\beta)
\qquad
\pi^\wedge_A\circ\pi^\wedge_W=
\pi_{E^\wedge}\circ\bomega^\wedge_0
$$
(where $\bomega^\wedge_0:W(E^\wedge)\to E^\wedge$ is
the $0$-th ghost component map, and
$\pi_{E^\wedge}:E^\wedge\to E^\wedge/a_0E^\wedge$ is the
projection) in light of which, it suffices to show that
$$
\pi_{E^\wedge}\circ\beta=\iota\circ\pi_A\circ\bar u_{A(E)}.
$$
The latter follows by a simple inspection of the
definition of $\beta$ : details left to the reader.
\end{pfclaim}
\end{proof}

\sset\subsubsection{}\label{subsec_equivalence}
For any prime $p$, consider the following categories :
\begin{itemize}
\item
The category $\cE$ whose objects are all the pairs
$(E,\cI)$ where $E$ is a perfect topological
$\F_p$-algebra and a P-ring, and $\cI\subset W(E)$ is a
distinguished ideal. The morphisms $(E,\cI)\to(E',\cI')$
are the continuous ring homomorphisms $f:E\to E'$ such that
$W(f)(\cI)\subset\cI'$.
\item
The category $\cA$ of all P-rings $A$ such that
$\Ker\,u_A$ is a distinguished ideal of $W(A)$. The
morphisms in $\cA$ are all the continuous ring
homomorphisms.
\item
The full subcategory $\cE_p$ of $\cE$ whose objects
are all the pairs $(E,\cI)$ such that the topology
of $E$ agrees with the $\bomega_0(\cI)$-adic topology,
and the full subcategory $\cA_p$ of $\cA$ consisting
of those objects $A$ whose topology agrees with the
$p$-adic topology.
\end{itemize}
Then proposition \ref{prop_special-case}(iii,iv)
yields an equivalence of categories
\set\begin{equation}\label{eq_equivalence}
\cE_p\isom\cA_p
\quad : \quad
(E,\cI)\mapsto A(E,\cI):=W(E)/\cI
\end{equation}
with quasi-inverse given by the rule :
$A\mapsto(\bE(A),\Ker\,u_A)$. In the following, we aim
to extend this equivalence to the whole of $\cE$. The
first step is :

\begin{proposition}\label{prop_bar-W-perfect}
With the notation of \eqref{subsec_equivalence}, let
$(E,\cI)$ be any object $\cE$. Then :
\begin{enumerate}
\item
$A(E,\cI)$ is a P-ring.
\item
Especially, $\cI$ is closed in the topology $\cT_{W(E)}$.
\end{enumerate}
\end{proposition}
\begin{proof} Clearly, it suffices to check (i). To this
aim, let $(a_n~|~n\in\N)\in\cI$ be any distinguished
element, $I\subset E$ be any ideal of definition,
$(x_0,\dots,x_r)$ a finite system of generators for $I$,
denote by $\cJ\subset W(E)$ the ideal generated by the
system $(\tau_E(x_0),\dots,\tau_E(x_r))$, and set
$\cJ':=\cJ+pW(E)$.
By proposition \ref{prop_morel}(ii), we know that the
topology of $A(E)$ agrees with the $\cJ'$-adic topology.
Now, by assumption $a_0$ is topologically nilpotent,
hence $a_0^k\in I$ for every sufficiently large $k\in\N$;
since the Frobenius endomorphism $\Phi_E$ is an automorphism
of the topological ring $E$, we may then replace
$I$ by $\Phi_E^{-t}I$ for some suitably large $t\in\N$,
after which we may assume that $a_0\in I^{(p)}$, and recall
that $I^{(p)}$ is generated by $(x_0^p,\dots,x_r^p)$.
Then, lemma \ref{lem_before-name}(ii) implies that 
\set\begin{equation}\label{eq_get-it-on-the-nose}
p\in\cJ^2A(E).
\end{equation}
and therefore $\cJ'A(E)=\!\cJ\!A(E)$; especially, the
topology of $A(E)$ agrees with the $\cJ$-adic topology.
We remark :

\begin{claim}\label{cl_con-game}
$\cI$ is a closed subset for the $\cJ$-adic topology on $W(E)$.
\end{claim}
\begin{pfclaim} By lemma \ref{lem_Witt-limit}(ii) and
\ref{lem_fontaine}, the ring $W(E)$ is complete and
separated for its $\cJ$-adic topology, so the assertion
follows from corollary \ref{cor_square-powers}.
\end{pfclaim}

Claim \ref{cl_con-game} implies that $A(E)$
is complete and separated for its $\cJ$-adic topology.
Moreover, the Frobenius endomorphism of $A(E)/\!\cJ\!A(E)$
is surjective, by virtue of \eqref{eq_mod-out-by-p};
taking into account \eqref{eq_get-it-on-the-nose}, the
proposition follows.
\end{proof}

To define a quasi-inverse for our extension of
\eqref{eq_equivalence}, we shall also need the following:

\begin{proposition}\label{prop_cond-for-perfectoid}
Let $A$ be any P-ring. We have :
\begin{enumerate}
\item
If\/ $\Ker\,u_A$ is the topological closure of a finitely
generated ideal, $\bE(A)$ is a P-ring.
\item
The following conditions are equivalent :
\begin{enumerate}
\item
$\Ker\,u_A$ is generated by any distinguished element
contained in it.
\item
$\Ker\,u_A$ is a principal ideal.
\item
$\Ker\,u_A$ is the topological closure of a principal ideal.
\end{enumerate}
\end{enumerate}
\end{proposition}
\begin{proof}(i): For any topological space $X$, and any
subset $S\subset X$, denote by $S^c$ the topological closure
of $S$ in $X$.
By assumption, there exists a finitely generated ideal
$\cI\subset\bA(A)$ such that $\cI^c=\Ker\,u_A$.
Let now $\bomega_0:\bA(A)\to\bE(A)$ be the $0$-th
ghost component; we notice :

\begin{claim}\label{cl_small-link}
(i)\ \
$\bomega_0(\Ker\,u_A)=\Ker\,(\bar u_{A/pA}:\bE(A)\to A/pA)$.
\begin{itemize}
\item[(ii)]
$(\Ker\,\bar u_{A/pA})^c=\Ker\,\bar u_{A/(pA)^c}$.
\end{itemize}
\end{claim}
\begin{pfclaim}(i): Say that $x\in\Ker\,\bar u_{A/pA}$,
and pick any $y\in\bA(A)$ such that $\bomega_0(y)=x$;
then $u_A(y)=pa$ for some $a\in A$. We pick $z\in\bA(A)$
such that $u_A(z)=a$; then $y-pz\in\Ker\,u_A$ and
$\bomega_0(y-pz)=x$, whence the contention.

(ii): The ideal $K:=(\Ker\,\bar u_{A/pA})^c$ contains
$\Ker\,\bar u_{A/pA}$, therefore
$$
K=\bar u_{A/pA}^{-1}(\bar u_{A/pA}(K))
$$
so $\bar u_{A/pA}(K)$ is a closed ideal of $A/pA$, by
lemma \ref{lem_was-third-cond}(i), hence
$(pA)^c/pA\subset\bar u_{A/pA}(K)$ and consequently
$K=\bar u^{-1}_{A/pA}((pA)^c/pA)$. On the other hand,
$\bar u_{A/(pA)^c}$ is the composition of $\bar u_{A/pA}$
with the projection $A/pA\to A/(pA)^c$, whence the claim.
\end{pfclaim}

Claims \ref{cl_image-and-closure} and \ref{cl_small-link}(i)
imply that $\bomega_0(\cI)^c=(\Ker\,\bar u_{A/pA})^c$. From claim
\ref{cl_small-link}(ii) we deduce that $\Ker\,\bar u_{A/(pA)^c}$
is the topological closure in $\bE(A)$ of a finitely generated
ideal. Next, let $I\subset A$ be any ideal of definition, and
$\cJ\subset\bE(A)$ any finitely generated ideal such
that $\bar u_{A/(pA)^c}(\cJ)=I/(pA)^c$; it follows easily
that $\cJ+\Ker\,\bar u_{A/(pA)^c}=\Ker\,\bar u_{A/I}$
(details left to the reader), and then also
$\Ker\,\bar u_{A/I}$ is the topological closure
of a finitely generated ideal of $\bE(A)$.
Lastly, we have a commutative diagram of continuous
and surjective ring homomorphisms
$$
{\diagram
\bE(A) \ar[r]^-{\bar u_{A/I}} \ar[d]_v &
A/I \ddouble \\
A/I \ar[r]^-{\Phi_{A/I}} & A/I
\enddiagram}
\qquad
\text{where $v:=\bar u_{A/I}\circ\Phi^{-1}_{\bE(A)}$}
$$
whence $v(\Ker\,\bar u_{A/I})=\Ker\,\Phi_{A/I}$.
Since the topology of $A/I$ is discrete, claim
\ref{cl_image-and-closure} implies that $\Ker\,\Phi_{A/I}$
is a finitely generated ideal, so the assertion
follows from proposition \ref{prop_third-equiv-cond}.

(ii): Clearly (a)$\Rightarrow$(b), and since $\Ker\,u_A$
is a closed subset of $W(A)$, we have as well
(b)$\Rightarrow$(c). Thus, it remains only to check
that (c)$\Rightarrow$(a). To this aim, let again
$I\subset A$ be any ideal of definition, denote
by $J\subset A$ the radical of $I$, and fix any
element $\underline b:=(b_n~|~n\in\N)\in\bA(A)$
such that $\Ker\,u_A$ is the topological closure of
$\underline b\bA(A)$ in $\bA(A)$. We notice that $A/J$
is a perfect $\F_p$-algebra, and its quotient topology
induced by the projection $\pi:A\to A/J$ is discrete.
We have a commutative diagram of continuous and surjective
ring homomorphisms :
$$
\xymatrix{ \bA(A) \ar[r]^-{u_A} \ar[d]_f & A \ar[d]^\pi \\
W(A/J) \ar[r]^-{\bomega_0} & A/J
}$$
whose bottom horizontal arrow is the $0$-th ghost map
(see \eqref{subsec_Witt-ghost}), and with $f:=W(\phi)$,
where $\phi:\bE(A)\to A/J$ is the composition of
$\bar u_{A/pA}:\bE(A)\to A/pA$ with the natural
projection $A/pA\to A/J$ (details left to the reader).

\begin{claim}\label{cl_missing-link}
$f(\Ker\,u_A)=pW(A/J)$.
\end{claim}
\begin{pfclaim} It is clear that $f(\Ker\,u_A)\subset pW(A/J)$,
so we need only show the converse inclusion. Now, according to
lemma \ref{lem_was-third-cond}(ii), we may find a distinguished
element $\underline a:=(a_n~|~n\in\N)$ in $\Ker\,u_A$. Since
$a_0$ is topologically nilpotent in $\bE(A)$, the element
$\bar u_{A/pA}(a_0)$ is topologically nilpotent in $A/pA$, and
therefore its image vanishes in $A/J$; in other words,
$a_0\in\Ker\,\phi$. However,
$\underline a=\tau_A(a_0)+p\cdot\underline u$, for some
$\underline u\in\bA(A)^\times$, so that
$$
f(\underline a)=f\circ\tau_A(a_0)+p\cdot f(\underline u)=
\tau_{A/J}\circ\phi(a_0)+p\cdot f(\underline u)=
p\cdot f(\underline u)
$$
(see \eqref{eq_Teich-functorial}) whence the claim.
\end{pfclaim}

In light of claims \ref{cl_image-and-closure} and
\ref{cl_missing-link} and example
\ref{ex_discrete-Witt}(i) we deduce that
$f(\underline b)\cdot W(A/J)$ is a dense subset of
$pW(A/J)$, for the $p$-adic topology on $W(A/J)$;
especially, we have
$$
pW(A/J)=p^2W(A/J)+f(\underline b)\cdot W(A/J).
$$
Then, since $p$ lies in the Jacobson radical of $W(A/J)$,
Nakayama's lemma shows that
$f(\underline b)\cdot W(A/J)=pW(A/J)$, which in turn
implies that $\phi(b_0)=0$ and $\phi(b_1)$ is invertible
in $A/J$.

\begin{claim}\label{cl_b-is-disting}
$\underline b$ is a distinguished element of $\bA(A)$.
\end{claim}
\begin{pfclaim} First, since $J$ is contained in the
Jacobson radical of $A$ (remark
\ref{rem_someth-on-bdd-in-Z-lin}(v)) and
$\phi(b_1)\in(A/J)^\times$, we deduce that
$\bar u_A(b_1)\in A^\times$, and therefore
$b_1\in\bE(A)^\times$, by remark \ref{rem_general}(iii).
Likewise, we have $\bar u_A(b_0)\in JA$, which
easily implies that $b_0$ is topologically nilpotent
in $\bE(A)$ (details left to the reader).
\end{pfclaim}

Now, from (i) and assumption (ii.c) we see that $\bE(A)$
is a P-ring; from claim \ref{cl_b-is-disting} and proposition
\ref{prop_bar-W-perfect}(ii) it then follows that
$\underline b\bA(A)$ is a closed ideal of $\bA(A)$, and
therefore it must coincide with $\Ker\,u_A$. Lastly, let
$\underline a$ be any other
distinguished element of $\bA(A)$ contained in $\Ker\,u_A$;
taking into account remark \ref{rem_distinguished}(i) and
again claim \ref{cl_b-is-disting}, we conclude that
$\underline a\bA(A)=\underline b\bA(A)$, whence (ii.a).
\end{proof}

\subsection{Perfectoid rings}\label{sec_now-perfectoid}
We are now ready to introduce our generalizations of Scholze's
perfectoid rings that will intervene in the proof of the
log regular version of almost purity.

\begin{definition} Let $A$ be any P-ring. We say
that $A$ is {\em perfectoid} if it fulfills any of
the three equivalent conditions of proposition
\ref{prop_cond-for-perfectoid}(ii).
\end{definition}

\begin{example}\label{ex_perfectoid}
(i)\ \
Let $(A,\cT)$ be any topological $\F_p$-algebra. Then
$A$ is perfectoid if and only if it is perfect, $\cT$
is complete and separated, and there exists a finitely
generated ideal $I\subset A$ such that $\cT$ agrees with
the $I$-adic topology. Indeed, under these assumptions,
$A$ is a P-ring (remark \ref{rem_why-only-now}), and the
map $\bar u_A:\bE(A)\to A$ is an isomorphism of topological
rings that identifies $u_A$ with the $0$-th ghost map
$\bomega_0:W(A)\to A$, whose kernel is the principal ideal
$pW(A)$, by \eqref{subsec_perfect-case}, so $A$ is perfectoid.
Conversely, clearly $p\in\Ker\,u_A$, and $p$ is obviously
a distinguished element of $\bA(A)$, so if $A$ is a
perfectoid $\F_p$-algebra we must have $\Ker\,u_A=pW(A)$;
moreover, $u_A$ is an open ring homomorphism (lemma
\ref{lem_was-third-cond}(i)), so it induces an isomorphism
$\bE(A)\isom A$ of topological rings, especially $(A,\cT)$
is a perfect topological $\F_p$-algebra, as claimed.

(ii)\ \
Let $E$ be any perfectoid $\F_p$-algebra, and
$\underline a:=(a_n~|~n\in\N)$ any distinguished element
of $W(E)$. Then $A(E):=W(E)/\underline a W(E)$ is
perfectoid for the quotient topology induced by $W(E)$.
Indeed, $A(E)$ is a P-ring, by proposition
\ref{prop_bar-W-perfect}(i); next, lemma \ref{lem_fontaine}
implies that $E$ is complete and separated for the
$a_0$-adic topology, and therefore proposition
\ref{prop_special-case}(ii,iv) applies, and shows that
$\Ker\,u_{A(E)}$ is a principal ideal, whence the contention.
Moreover, in this case the perfectoid ring $E$ can be
recovered from $A(E)$ : indeed, we have a natural
isomorphism of topological rings
$$
\bE(A(E))\isom E
$$
obtained as follows. First, since the $0$-th ghost
component $W(E)\to E$ is an open and surjective map, the
same holds for the projection
$\alpha:A(E)\to A(E)/pA(E)\isom E/a_0E$, where $E/a_0E$ is
endowed with the quotient topology induced by $E$.
Then, $\alpha$ induces an isomorphism
$\bE(\alpha):\bE(A(E))\isom\bE(E/a_0E)$ of topological
rings, by \eqref{sec_def-A-tilde}. But since $a_0$
is topologically nilpotent, the projection
$\beta:\bE(E)\to\bE(E/a_0E)$ is an isomorphism of
topological rings as well (theorem \ref{th_fontaine}),
and the same holds for the map $\bar u_E:\bE(E)\to E$,
since $E$ is perfect. Thus, the sought isomorphism is
the composition $\bar u_E\circ\beta^{-1}\circ\bE(\alpha)$.
\end{example}

Let us also point out the following complement to
proposition \ref{prop_third-equiv-cond} :

\begin{corollary}\label{cor_jackob}
Let $A$ be any perfectoid ring. Then $\bar\Phi_{A/J}$
is an isomorphism for every special ideal of definition
$J$ of $A$.
\end{corollary}
\begin{proof} Set $E:=\bE(A)$; we may assume that
there exists a distinguished element
$\underline a:=(a_n~|~n\in\N)\in W(E)$ such that $A$ is
the ring $A(E)$ as defined in \eqref{subsec_several-tops},
and $u_A:W(E)\to A(E)$ is the projection. Moreover, we may
find a finite system $\underline x{}_\bullet:=
(\underline x{}_0,\dots,\underline x{}_r)$ of elements
of\/ $W(E)$, such that $J$ is the image in $A(E)$ of
the ideal generated by $\underline x{}_\bullet$.
We may write
\set\begin{equation}\label{eq_pirana}
\underline x{}_i=\tau_\bE(t_i)+p\cdot\underline y{}_i
\qquad
\text{where $t_i\in E$ and $\underline y{}_i\in W(E)$
for every $i=0,\dots,r$}.
\end{equation}
In view of lemma \ref{lem_basic-cong}(i), we deduce
that
\set\begin{equation}\label{eq_taus-and-ts}
\underline x{}_i^p\equiv\tau_E(t_i)^p\pmod{p^2W(E)}
\qquad
\text{for $i=0,\dots,r$}.
\end{equation}
Let $\cJ\subset W(E)$ be the ideal generated by
$\tau_E(t_0),\dots,\tau_E(t_r)$; we notice :

\begin{claim}\label{cl_equal-in-bar}
$p\in\cJ^{(p)}A(E)$ and $\cJ A(E)=J$.
\end{claim}
\begin{pfclaim} It follows easily from \eqref{eq_taus-and-ts}
that there exists $\underline z\in W(E)$ with
$$
p\cdot u_A(1+p\cdot\underline z)\in\cJ^{(p)}A(E)
$$
(recall that $\cJ^{(p)}$ is the ideal generated by
$\tau_E(t_0)^p,\dots,\tau_E(t_r)^p$ : see definition
\ref{def_perfectoid}). However, $p$ is topologically
nilpotent in $W(E)$, therefore
$1+p\cdot\underline z\in W(E)^\times$, so
$p\in\cJ^{(p)}A(E)$. The other assertion follows
directly from the first one and \eqref{eq_pirana} : details
left to the reader.
\end{pfclaim}

Let $I\subset E$ be the ideal generated by $t_0,\dots,t_r$,
and set $E_0:=E/a_0E$; taking into account lemma
\ref{lem_before-name}(ii) and claim \ref{cl_equal-in-bar},
we deduce that $a_0\in I^{(p)}$. Also, it is clear that the
natural isomorphism \eqref{eq_mod-out-by-p} maps the image
of $J$ (resp. of $J^{(p)}$) onto $IE_0$ (resp. $I^{(p)}E_0$).
But it also follows that $E_0/IE_0=E/I$ and
$E_0/I^{(p)}E_0=E/I^{(p)}$; lastly, since $E$ is perfect, the
Frobenius endomorphism $\Phi_E$ maps $I$ onto $I^{(p)}$, whence
the contention.
\end{proof}

\begin{remark}\label{rem_nice-topology}
(i)\ \
We may now say that the category $\cA$ of
\eqref{subsec_equivalence} is the {\em category of
perfectoid rings}, and we shall denote it henceforth by
$$
\Perf.
$$
More generally, if $A$ is any perfectoid ring, we shall
write
$$
A\tdu\Perf
$$
for the category $A/\Perf$ of {\em perfectoid $A$-algebras},
whose objects are the continuous ring homomorphisms $A\to B$
with $B$ perfectoid. The objects of the category $\cE$ of
\eqref{subsec_equivalence} are the pairs $(E,\cI)$ consisting
of a perfectoid $\F_p$-algebra and a distinguished ideal
of $W(E)$. Taking into account example \ref{ex_perfectoid},
we see that  \eqref{subsec_equivalence} generalizes
{\em verbatim} to an equivalence
$$
\cE\isom\Perf
\quad : \quad
(E,\cI)\mapsto A(E,\cI):=W(E)/\cI
$$
with a natural quasi-inverse given by the functor $\bE$,
which restricts to equivalences
$$
\bE:A\tdu\Perf\isom\bE(A)\tdu\Perf
\qquad
\text{for every perfectoid ring $A$}.
$$
Hereafter, we shall study the stability of the class of
perfectoid rings under some standard operations. The
following observations (iii) and (iv) show stability
under completion with respect to the $J$-adic topology
corresponding to an ideal $J\subset A$ of finite type
with $pA\subset J$.

(ii)\ \
Let $A$ be any perfectoid ring, set $\bE:=\bE(A)$, and
let $\alpha_\bullet:=(\alpha_n~|~n\in\N)$ be any distinguished
element in $\Ker\,u_A$. Notice that $\bar u_{A/pA}:\bE\to A/pA$
induces an isomorphism of topological rings :
$$
\omega:\bE/\alpha_0\bE\isom A/pA
$$
for the quotient topologies induced by the projections
$A\to A/pA$ and $\bE\to\bE/\alpha_0$. Indeed, $\bar u_{A/pA}$
is open and surjective (lemma \ref{lem_was-third-cond}(i)),
and $\alpha_\bullet$ generates $\Ker\,u_A$, so $\alpha_0$
generates the kernel of $u_{A/pA}$, by virtue of lemma
\ref{lem_drop-conditions}(i), whence the contention.

(iii)\ \
The isomorphism $\omega$ of (ii) induces a natural
bijection :
$$
\{\text{ideals $\cJ\subset\bE$}~|~\alpha_0\bE\subset\cJ\}
\leftrightarrow
\{\text{ideals $J\subset A$}~|~pA\subset J\}.
$$
Namely, $\cJ\subset\bE$ and $J\subset A$ correspond
under this bijection, if and only if
$\omega(\cJ/\alpha\bE)=J/pA$. Then, clearly $J$ is
finitely generated if and only if the same holds for
$\cJ$. Let $\cT_J$ (resp. $\cT_{\!\!\cJ}$) denote the
$J$-adic (resp. $\cJ$-adic) topology on $A$ (resp. on $\bE$);
taking into account remark \ref{rem_topology-of-E}(ii),
we deduce that if $J$ is finitely generated, then the
topology of $\bE(A,\cT_J)$ agrees with $\cT_{\!\!\cJ}$.
Combining with example \ref{ex_perfectoid}(i),
lemma \ref{lem_still-c-adic}(iv), example
\ref{ex_discrete-Witt}(ii) and remark \ref{rem_topology-of-E}(v),
we conclude that the completion
$\bE(A,\cT_J)^\wedge$ of $\bE(A,\cT_J)$ is perfectoid.

(iv)\ \
In the situation of (iii), suppose that $J$ and $\cJ$ are
finitely generated; notice that, since $\alpha_0\in\cJ$,
the image of $\alpha_\bullet$ is distinguished in
$W(\bE(A,\cT_J)^\wedge)$, and let
$\cI:=\alpha_\bullet W(\bE(A,\cT_J)^\wedge)$. We claim that
there is a natural isomorphism of topological rings
$$
W(\bE(A,\cT_J)^\wedge)/\cI
\isom(A_J,\cT_J)^\wedge
$$
for the quotient topology on $W(\bE(A,\cT_J)^\wedge)/\cI$
induced by $W(\bE(A,\cT_J)^\wedge)$. Especially, this shows
that $(A_J,\cT_J)^\wedge$ is perfectoid. For the proof, fix
a finite system $(\beta_1,\dots,\beta_k)$ of generators of
$\cJ$, and let $\cJ_W\subset W(\bE)$ be the ideal generated by
the system $(p,\tau_\bE(\beta_1),\dots,\tau_\bE(\beta_k))$;
by proposition \ref{prop_morel}(ii), the topology
$\cT_W$ of $W(\bE,\cT_{\!\!\cJ})$ agrees with the $\cJ_W$-adic
topology. By construction we have $u_A(\cJ_W)/pA=J/pA$, 
and consequently $u_A(\cJ_W)=J$, so the quotient topology
induced by $\cT_W$ on $A$ via $u_A$ agrees with $\cT_J$. Lastly,
in view of (iii) and proposition \ref{prop_bar-W-perfect}(ii)
we know that $\cI$ is a closed ideal of $W(\bE(A,\cT_J)^\wedge)$.
Taking into account proposition \ref{prop_replaces-Mat-Th.8.1}(v)
and lemma \ref{lem_Witt-limit}(iv), the assertion follows.

We can summarize these observations in the following :
\end{remark}

\begin{proposition}\label{prop_change-topol}
In the situation of proposition {\em\ref{prop_P-rings}},
the following holds :
\begin{enumerate}
\item
If $(A,\cT_I)$ is perfectoid, the same holds for
$(A,\cT_p)$.
\item
If $(A,\cT_p)$ is perfectoid, the same holds for
the completion $(A_I,\cT_I)^\wedge$ of $(A,\cT_I)$.
\end{enumerate}
\end{proposition}
\begin{proof}(i) follows immediately from proposition
\ref{prop_P-rings}(i).

(ii): It suffices to apply remark \ref{rem_nice-topology}(iv)
to the ideal $I\subset A$ and the corresponding ideal
$\cJ:=\bar u_{A/pA}^{-1}(I/pA)$ of $\bE(A)$.
\end{proof}

The next observation establishes the stability of $\Perf$
under completed tensor products :

\begin{proposition}\label{prop_stabil-cplted-tensors}
Let $A_0$ be a perfectoid ring and $A_1$, $A_2$ two
perfectoid $A_0$-algebras. Then:
\begin{enumerate}
\item
The topological ring $A_3:=A_1\,\hat\otimes_{A_0}A_2$ is
perfectoid.
\item
There exists a natural isomorphism
$\bE(A_3)\isom\bE(A_1)\,\hat\otimes_{\bE(A_0)}\bE(A_2)$
in $\Perf$.
\item
Especially, all finite coproducts are representable in the
category $\Perf$.
\end{enumerate}
\end{proposition}
\begin{proof}(i): Set $\bE_i:=\bE(A_i)$ for $i=0,1,2$, and
pick a distinguished element $\underline\alpha\in W(\bE_0)$
that generates the kernel of $u_{A_0}$; then the image of
$\underline\alpha$ in $W(\bE_i)$ is still distinguished,
and we know that $\underline\alpha W(\bE_i)$ is the kernel
of $u_{A_i}$ also for $i=1,2$. Moreover, it is clear from
remark \ref{rem_tensor-Witt}(iii) that the $\F_p$-algebra
$\bE_3:=\bE_1\,\hat\otimes_{\bE_0}\bE_2$ is perfectoid,
and then the same holds for the topological ring
$A(\bE_3):=W(\bE_3)/\underline\alpha W(\bE_3)$, by
virtue of example \ref{ex_perfectoid}(ii). Lastly,
the isomorphism \eqref{eq_tensor-Witt} induces an
isomorphism of topological rings
$$
A(\bE_3)\isom A_1\,\hat\otimes_{A_0}A_2
$$
whence the assertion.

(ii): Example \ref{ex_perfectoid}(ii) also identifies
naturally $\bE_3$ with $\bE(A(\bE_3))$, whence the contention.

(iii) is clear, since complete tensor products represent
these coproducts, by \eqref{subsec_tensor-topol}.
\end{proof}

\sset\subsubsection{}\label{subsec_double-completion}
For any ring $A$ and any ideal $I\subset A$, let
$\cT_I$ be the $I$-adic topology on $A$, and denote
$(A^\wedge_I,\cT^\wedge_I)$ the separated completion
of $(A,\cT_I)$. Notice that if $I'\subset A$ is another
ideal with $I'\subset I$, the identity map of $A$ yields
a continuous map $(A,\cT_{I'})\to(A,\cT_I)$, whose
completion is a natural continuous ring homomorphism
\set\begin{equation}\label{eq_I-and-I-prime}
(A^\wedge_{I'},\cT^\wedge_{I'})\to(A^\wedge_I,\cT^\wedge_I).
\end{equation}
Now, suppose that $I,J\subset A$ are any two ideals;
we deduce a commutative diagram
\set\begin{equation}\label{eq_old-acquaintance}
{\diagram
(A^\wedge_{I\cap J},\cT^\wedge_{I\cap J})
\ar[r]^-{\phi_I} \ar[d]_{\phi_J} &
(A^\wedge_I,\cT^\wedge_I) \ar[d]^{\beta_I} \\
(A^\wedge_J,\cT^\wedge_J) \ar[r]^-{\beta_J} &
(A^\wedge_{I+J},\cT^\wedge_{I+J})
\enddiagram}
\end{equation}
whose arrows are the continuous ring homomorphisms
\eqref{eq_I-and-I-prime}. Moreover, for every $k\in\N$
let $A^\wedge_{I,k}:=A^\wedge_I/J^kA^\wedge_I$, and endow
$A^\wedge_{I,k}$ with the quotient topology $\cT_{I,k}$
induced by $\cT^\wedge_I$ via the natural projection
$A^\wedge_I\to A^\wedge_{I,k}$. We set
\set\begin{equation}\label{eq_second-completion}
(A^\wedge_{I,J},\cT^\wedge_{I,J}):=
\lim_{k\in\N}\,(A^\wedge_{I,k},\cT_{I,k}).
\end{equation}
For every $k\in\N$, let $\cT_{d,k}$ be the discrete
topology on $A/(I+J)^k$;  the system of continuous
projections
$((A^\wedge_I,\cT^\wedge_I)\to(A/(I+J)^k,\cT_d)~|~k\in\N)$
factors uniquely through a system of continuous ring
homomorphisms
$((A^\wedge_{I,k},\cT_{I,k})\to(A/(I+J)^k,\cT_d)~|~k\in\N)$,
whose limit is a natural continuous ring homomorphism
\set\begin{equation}\label{eq_double-completion}
(A^\wedge_{I,J},\cT^\wedge_{I,J})\to
(A^\wedge_{I+J},\cT^\wedge_{I+J}).
\end{equation}

\begin{lemma}\label{lem_double-completion}
In the situation of \eqref{subsec_double-completion},
we have :
\begin{enumerate}
\item
\eqref{eq_double-completion} is an open and surjective map.
\item
Suppose that $A$ is either a noetherian ring or a perfect\/
$\F_p$-algebra, and that both $I$ and $J$ are finitely
generated. Then :
\begin{enumerate}
\item
\eqref{eq_old-acquaintance} is a cartesian diagram
of topological rings and a cocartesian diagram of
topological $A$-modules.
\item
\eqref{eq_double-completion} is an isomorphism of
topological rings.
\end{enumerate}
\end{enumerate}
\end{lemma}
\begin{proof} For every $n\in\N$, we have a natural
diagram of rings
\set\begin{equation}\label{eq_like-old-times}
{\diagram
A/(I^n\cap J^n) \ar[r] \ar[d] & A/I^n \ar[d] \\
A/J^n \ar[r] & A/(I^n+J^n)
\enddiagram}
\end{equation}
whence a complex of $A$-modules
$$
\Sigma_n
\qquad :\qquad
0\to A/(I^n\cap J^n)\to A/I^n\oplus A/J^n\to A/(I^n+J^n)\to 0
$$
and it is easily seen that $\Sigma_n$ is exact for
every $n\in\N$. Therefore, \eqref{eq_like-old-times}
is a cartesian diagram of discrete topological rings
(details left to the reader); let $\cT'_{d,n}$ be the
discrete topology on $A/(I^n\cap J^n)$ for every
$n\in\N$, and set
$$
(A',\cT'):=\lim_{n\in\N}\,(A/(I^n\cap J^n),\cT'_{d,n}).
$$
Notice as well that $(I+J)^{2n-1}\subset I^n+J^n\subset(I+J)^n$
for every $n\in\N$, so the $(I+J)$-adic topology on $A$
agrees with the linear topology defined by the descending
system of ideals $(I^n+J^n~|~n\in\N)$.
In light of example \ref{ex_lim_interchange}(ii), we
deduce that the limit of the system of these diagrams
is a cartesian diagram of topological rings
\set\begin{equation}\label{eq_same-again}
{\diagram
(A',\cT') \ar[r]^-{\phi'_I} \ar[d]_{\phi'_J} &
(A^\wedge_I,\cT^\wedge_I) \ar[d]^{\beta'_I} \\
(A^\wedge_J,\cT^\wedge_J) \ar[r]^-{\beta'_J} &
(A^\wedge_{I+J},\cT^\wedge_{I+J}).
\enddiagram}
\end{equation}
Moreover, since the inverse system $(A/(I^n\cap J^n)~|~n\in\N)$
has surjective transition maps, the limit of the system
of exact sequences $(\Sigma_n~|~n\in\N)$ is still exact
(see \cite[Lemma 3.5.3]{We}), and therefore
\eqref{eq_same-again} is also a cocartesian diagram of
$A$-modules. Furthermore, by proposition
\ref{prop_replaces-Mat-Th.8.1}(i,v) we have an induced
morphism of exact complexes
$$
\xymatrix{
\Sigma'_n \ar[d]_{\sigma_n} &
0 \ar[r] & (I^n\cap J^n)^\wedge \ar[r] \ar[d] &
I^n{}^\wedge\oplus J^n{}^\wedge \ar[r] \ar[d] &
(I^n+J^n)^\wedge \ar[r] \ar[d] & 0 \\
\Sigma &
0 \ar[r] & A' \ar[r] &
A^\wedge_I\oplus A^\wedge_J \ar[r]^-{\beta'_{I,J}} &
A^\wedge_{I+J} \ar[r] & 0
}$$
where :
\begin{itemize}
\item
$(I^n\cap J^n)^\wedge$ is the topological closure
of $I^n\cap J^n$ in $(A',\cT')$, and likewise for
$I^n{}^\wedge\oplus J^n{}^\wedge$ and $(I^n+J^n)^\wedge$.
\item
$\Coker\,\sigma_n$ is naturally isomorphic
to $\Sigma_n$, for every $n\in\N$.
\item
$\beta'_{I,J}$ is the sum of $\beta'_I$ and $\beta'_J$.
\end{itemize}
Especially, $\beta'_{I,J}$ is a continuous, open and
surjective map, for the product topology
$\cT^\wedge_I\times\cT^\wedge_J$ on
$A^\wedge_I\oplus A^\wedge_J$, and it follows easily
that \eqref{eq_same-again} is a cocartesian diagram
of topological $A$-modules.

(i): The system of natural maps
$(A^\wedge_I\to A^\wedge_{I,k}\leftarrow A/J^k~|~k\in\N)$
yields ring homomorphisms
$$
A^\wedge_I\xrightarrow{\ \psi_I\ }(A^\wedge_I)^\wedge_J
\xleftarrow{\ \psi_J\ } A^\wedge_J.
$$
There follows a map of abelian groups
$\gamma:A^\wedge_I\oplus A^\wedge_J\to (A^\wedge_I)^\wedge_J$
which is continuous for the topology
$\cT^\wedge_I\times\cT^\wedge_J$ and a simple inspection
shows that $\eqref{eq_double-completion}\circ\gamma$
equals $\beta'_{I,J}$. Now, let $K\subset(A^\wedge_I)^\wedge_J$
be any open ideal; then $K':=\gamma^{-1}K$ is an open
subgroup of $A^\wedge_I\oplus A^\wedge_J$ and therefore
$\beta'_{I,J}K'$ is an open subgroup of $A^\wedge_{I+J}$
contained in the image $K''$ of $K$ under
\eqref{eq_double-completion}, so $K''$ is open as well,
whence (i).

(ii.a): The system of projections
$(A/(I\cap J)^n\to A/(I^n\cap J^n)~|~n\in\N)$ yields
a natural continuous ring homomorphism
$$
\nu:(A^\wedge_{I\cap J},\cT^\wedge_{I\cap J})\to
(A',\cT')
$$
and a simple inspection shows that $\phi'_I\circ\nu=\phi_I$
and $\phi'_J\circ\nu=\phi_J$. To conclude the proof of (ii.a),
it then suffices to check :

\begin{claim} Under the assumptions of (ii), the map
$\nu$ is an isomorphism of topological rings.
\end{claim}
\begin{pfclaim} Since
$(IJ)^{2n}\subset(I\cap J)^{2n}\subset(IJ)^n$ for every
$n\in\N$, the $(I\cap J)$-adic and $IJ$-adic topologies
coincide on $A$. On the other hand, in case $A$ is
noetherian, the Artin-Rees lemma (\cite[Th.8.5]{Mat})
implies that for every $n\in\N$ there exists $m\in\N$
such that $I^m\cap J^n\subset I^nJ^n\subset(I\cap J)^n$,
from which the claim follows easily.

Thus, suppose $A$ is perfect and both $I$ and $J$
are finitely generated, so the same holds for $IJ$,
and therefore lemma \ref{lem_mon-fract-powers}(iv) says
that the $IJ$-adic topology on $A$ agrees with the
linear topology defined by the cofiltered system of ideals
$((IJ)^{\La n\Ra}A~|~n\in\N)$. Likewise, the linear topology
on $A$ defined by the cofiltered system of ideals
$(I^n\cap J^n~|~n\in\N)$ agrees with the one defined by
the cofiltered system of ideals
$(I^{\La n\Ra}A\cap J^{\La n\Ra}A~|~n\in\N)$. Hence, it
suffices to prove that
$$
I^{\La pn\Ra}A\cap J^{\La pn\Ra}A\subset (IJ)^{\La n\Ra}A
\qquad
\text{for every $n\in\N$}.
$$
However, if $x\in I^{\La pn\Ra}A\cap J^{\La pn\Ra}A$, we may
write $x=x^{1/p}\cdot x^{(p-1)/p}\in I^{\La n\Ra}J^{\La(p-1)n/p\Ra}A
\subset(IJ)^{\La n\Ra}A$ whence the contention.
\end{pfclaim}

(ii.b): For every $h,k\in\N$, let $A_{h,k}:=A/(I^h+J^k)$,
and endow $A_{h,k}$ with the discrete topology $\cT_{h,k}$;
moreover, set
$$
(C_k,\cT_{C,k}):=\lim_{h\in\N}\,(A_{h,k},\cT_{h,k})
\qquad
\text{for every $k\in\N$}
$$
(where the transition maps $A_{h+1,k}\to A_{h,k}$ are the
natural projections). Since $I$ and $J$ are both finitely
generated, for every $h,k\in\N$ there exists $n\in\N$ such
that $(I+J)^n\subset I^h+J^k$, so the natural map
\set\begin{equation}\label{eq_third-completion}
(A^\wedge_{I+J},\cT^\wedge_{I+J})\to
\lim_{h,k\in\N}\,(A_{h,k},\cT_{h,k})\isom
\lim_{k\in\N}\,(C_k,\cT_{C,k})
\end{equation}
is an isomorphism of topological rings (see example
\ref{ex_lim_interchange}(ii)). On the other hand, the
projection $A\to A/J^k$ extends uniquely to a continuous
map $\phi_k:(A^\wedge_I,\cT^\wedge_I)\to(C_k,\cT_{C,k})$,
and the latter factors uniquely through a continuous
ring homomorphism
$$
\bar\phi_k:(A^\wedge_{I,k},\cT_{I,k})\to(C_k,\cT_{C,k})
\qquad
\text{for every $\in\N$}.
$$
With this notation, a simple inspection shows that
\eqref{eq_second-completion} and \eqref{eq_third-completion}
identify \eqref{eq_double-completion} with the limit
of the system of maps $(\bar\phi_k~|~k\in\N)$. However,
proposition \ref{prop_replaces-Mat-Th.8.1}(i,v) says that
$\phi_k$ is surjective and its kernel is the topological
closure $(J^kA^\wedge_I)^c$ of
$J^kA^\wedge_I$ in $(A^\wedge_I,\cT^\wedge_I)$, for every
$k\in\N$.
Thus, in order to show the lemma, it suffices to check
that the linear topology on $A^\wedge_I$ defined by the
system of ideals $((J^kA^\wedge_I)^c~|~k\in\N)$ agrees
with the $JA^\wedge_I$-adic topology. In case $A$ is
noetherian, this is clear, since in that case $J^kA^\wedge_I$
is already closed in $A^\wedge_I$ (\cite[Th.8.11]{Mat}).
For the case where $A$ is perfect, we remark :

\begin{claim}\label{cl_B-is-A-cplt}
Let $(B,\cT_B)$ be a perfect, complete and separated
topological $\F_p$-algebra, $I,J\subset B$ two ideals
of finite type, such that $\cT_B$ agrees with the
$I$-adic topology. For every $\lambda\in\N[1/p]$
let $J^{\La\lambda\Ra}B^c$ be the topological closure of
$J^{\La\lambda\Ra}B$ in $B$ (notation of
\eqref{subsec_mon-fract-powers}). Then
$$
J^{\La\lambda\Ra}B^c\subset
\bigcap_{n\in\N}(J+I^n)^{\La\lambda\Ra}B\subset
J^{\La\lambda'\Ra}B
\qquad
\text{for every $\lambda'<\lambda$}.
$$
\end{claim}
\begin{pfclaim} We show first that
$J^{\La\lambda\Ra c}B\subset J^{\La\lambda'\Ra}B$. Indeed,
pick a finite system of generators $(b_1,\dots,b_k)$
(resp. $(b_{k+1},\dots,b_s)$) for $J$ (resp. for $I$),
and let $x\in J^{\La\lambda\Ra}B^c$ be any element; we can
write 
$$
x=\sum_{n\in\N}x_n
\qquad\text{with}\qquad
x_n\in J^{\La\lambda\Ra}B\cap I^{\La n\Ra}B
\qquad
\text{for every $n\in\N$}.
$$
Choose a strictly positive $\eps\in\N[1/p]$ such that
$\lambda'':=(1-\eps)\cdot\lambda>\lambda'$, and notice that
$$
x_n=x_n^{1-\eps}\cdot x_n^\eps\in
J^{\La\lambda''\Ra}\cdot I^{\La n\eps\Ra}B
\qquad
\text{for every $n\in\N$}.
$$
By definition, for every $n\in\N$ there exist :
\begin{itemize}
\item
a finite set $S_n\subset\N[1/p]^{\oplus k}$ with
$\mu_1+\cdots+\mu_k=\lambda''$ for every
$\mu:=(\mu_1,\dots,\mu_k)\in S_n$
\item
a system $(a_\mu~|~\mu\in S_n)$ of elements of $B$ such that
$$
x_n^{1-\eps}=\sum_{\mu\in S_n}a_\mu b^\mu
\qquad
\text{where $b^\mu:=b_1^{\mu_1}\cdots b_k^{\mu_k}$ for every
$\mu\in S_n$}.
$$
\end{itemize}
Now, choose $N\in\N$ such that $\lambda''-\lambda'\geq kp^{-N}$
and define $\bar\mu$ and $\mu^*$ as in the proof of
proposition \ref{prop_morel}, so that
$\bar\mu\in S:=
\{\nu\in p^{-N}\N^{\oplus k}~|~\lambda''\geq\nu_1+\cdots+\nu_k>\lambda'\}$
for every $n\in\N$ and every $\mu\in S_n$. It follows that
$$
x=\sum_{\nu\in S}b^\nu c_\nu
\qquad\text{where}\qquad
c_\nu:=
\sum_{n\in\N}x_n^\eps\cdot\sum_{\substack{ \mu\in S_n \\ 
                                     \bar\mu=\nu}}a_\mu b^{\mu^*}.
$$
Clearly $b^\nu\in J^{\La\lambda'\Ra}B$, and by lemma
\ref{lem_mon-fract-powers}(iv) the series $c_\nu$ converges
in the $I$-adic topology of $B$ for every $\nu\in S$;
also it is easily seen that $S$ is a finite set, whence the
contention.

Next, for any $n\in\N$, the ideal $(J+I^n)^{\La\lambda\Ra}B$
is generated by all monomials of the form
$$
b^\mu:=b_1^{\mu_1}\cdots b_s^{\mu_s}
\qquad\text{such that}\qquad
\lambda_1+n^{-1}\cdot\lambda_2=\lambda
\quad\text{where}\quad
\lambda_1:=\sum_{i=1}^k\mu_i
\quad
\lambda_2:=\sum_{i=k+1}^s\mu_i
$$
and where $\mu:=(\mu_1,\dots,\mu_s)$ is any sequence of
elements of $\N[1/p]$. Fix also $\lambda''\in\N[1/p]$
with $\lambda'<\lambda''<\lambda$. Now, for any such
$\mu$, we have either $\lambda_1\geq\lambda''$, in which
case $b^\mu\in J^{\La\lambda''\Ra}B$, or else
$\lambda_2>n\cdot(\lambda-\lambda'')$. We conclude that
$$
\bigcap_{n\in\N}(J+I^n)^{\La\lambda\Ra}B\subset
\bigcap_{n\in\N}
(J^{\La\lambda''\Ra}B+I^{\La n\cdot(\lambda-\lambda'')\Ra}B)=
J^{\La\lambda''\Ra}B^c
$$
where the last identity follows from lemma
\ref{lem_mon-fract-powers}(iv). But we know already that
$J^{\La\lambda''\Ra}B^c\subset J^{\La\lambda'\Ra}B$, whence
the claim.
\end{pfclaim}

To conclude, it suffices now to apply claim
\ref{cl_B-is-A-cplt} with $B:=A^\wedge_I$ and
invoke lemma \ref{lem_mon-fract-powers}(iv) (details
left to the reader).
\end{proof}

\begin{theorem}\label{th_double-completion}
In the situation of \eqref{subsec_double-completion},
let $\cT_p$ be the $p$-adic topology of $A$, and suppose:
\begin{enumerate}
\alphaenu
\item
$(A,\cT_p)$ is perfectoid.
\item
$p^N\in I\cap J$ for every sufficiently large $N\in\N$.
\item
$I$ and $J$ are finitely generated ideals of $A$.
\end{enumerate}
Then we have :
\begin{enumerate}
\item
\eqref{eq_old-acquaintance} is a cartesian diagram
of rings and a cocartesian diagram of $A$-modules.
\item
\eqref{eq_double-completion} is an isomorphism of
topological rings.
\end{enumerate}
\end{theorem}
\begin{proof}(i): Let $\pi\in A$ be as in lemma
\ref{lem_perfectoid}(iv); quite generally, if $K\subset A$
is any finitely generated ideal containing $p^N$ (for some
$N\in\N$), then it is easily seen that the $K$-adic topology
on $A$ agrees with the $(K+\pi A)$-adic topology. In view of
our assumptions (b) and (c) we may therefore replace $I$ and
$J$ by $I+\pi A$ and respectively $J+\pi A$, and assume
from start that $\pi\in I\cap J$. Fix a finite system 
$(\bar b_1,\dots,\bar b_n)$ (resp. $(\bar b_{n+1},\dots,\bar b_m)$)
of generators of the ideal $I/pA$ (resp. of $J/pA$) of $A/pA$,
and let also $\underline\alpha:=(\alpha_n~|~n\in\N)$ be a
distinguished element of $\bA:=\bA(A,\cT_p)$ that generates
$\Ker\,u_A$. Pick elements $\beta_1,\dots,\beta_m\in\bE:=\bE(A)$
such that $\bar u_{A/pA}(\beta_i)=\bar b_i$ for $i=1,\dots,m$,
denote by $\cI_\bE$ (resp. $\cJ_\bE$) the ideal of
$\bE$ generated by the system $(\alpha_0,\beta_1,\dots,\beta_n)$
(resp. $(\alpha_0,\beta_{n+1},\dots,\beta_m)$), and
let $(\bE^\wedge_\cI,\cT^\wedge_\cI)$ (resp.
$(\bE^\wedge_\cJ,\cT^\wedge_\cJ)$, resp.
$(\bE^\wedge_{\cI\cap\cJ},\cT^\wedge_{\cI\cap\cJ})$, resp.
$(\bE^\wedge_{\cI+\cJ},\cT^\wedge_{\cI+\cJ})$) be the
$\cI_\bE$-adic (resp. $\cJ_\bE$-adic, resp. $(\cI\cap\cJ)$-adic,
resp. $(\cI+\cJ)$-adic) completion of $\bE$. By lemma
\ref{lem_double-completion}(ii), we get a cartesian diagram
of topological rings
$$
\cE\qquad :\qquad
{\diagram (\bE^\wedge_{\cI\cap\cJ},\cT^\wedge_{\cI\cap\cJ})
\ar[r] \ar[d] & (\bE^\wedge_\cI,\cT^\wedge_\cI) \ar[d] \\
(\bE^\wedge_{\cJ},\cT^\wedge_{\cJ}) \ar[r] &
(\bE^\wedge_{\cI+\cJ},\cT^\wedge_{\cI+\cJ})
\enddiagram}$$
which is also cocartesian as a diagram of topological
abelian groups.

\begin{claim}\label{cl_still-cart}
$W(\cE)$ is still a cartesian diagram of topological rings,
and a cocartesian diagram of $\bA$-modules.
\end{claim}
\begin{pfclaim} The cartesian property follows from
remark \ref{rem_Witt-limit}(ii). To prove the cocartesian
property, it suffices to show that the map of abelian
groups deduced from $W(\cE)$
$$
\tau:W(\bE^\wedge_\cI)\oplus W(\bE^\wedge_\cJ)\to
W(\bE^\wedge_{\cI+\cJ})
$$
is surjective. Endow both target and source of $\tau$
with their $p$-adic filtrations, so that $\tau$ becomes
a map of filtered abelian groups, and for every $k\in\N$,
denote by $\gr^k\tau$ the map induced by $\tau$ on the
respective $k$-graded subquotients. Since both target
and source of $\tau$ are $p$-adically complete and
separated (proposition \ref{prop_Witt-is-complete}(iii)),
it then suffices to check that $\gr^k\tau$ is surjective
for every $k\in\N$ (\cite[Ch.III, \S2, n.8, Cor.2]{BouAC}).
Since multiplication by $p$ is an injective endomorphism
for both target and source of $\tau$, we are then further
reduced to the case where $k=0$, in which case $\gr^0\tau$
is naturally identified with the map
$\bE^\wedge_\cI\oplus\bE^\wedge_\cJ\to\bE^\wedge_{\cI+\cJ}$
deduced from $\cE$. But the latter map is indeed surjective,
since $\cE$ is cocartesian.
\end{pfclaim}

Now, remark \ref{rem_nice-topology}(iii) says that the
topology of $\bE(A,\cT_I)$ agrees with the $\cI_\bE$-adic
topology, so the completion $\bE_I^\wedge$ of $\bE(A,\cT_I)$
is isomorphic to $(\bE^\wedge_\cI,\cT^\wedge_\cI)$.
Likewise, the completion $\bE_J^\wedge$ of $\bE(A,\cT_J)$ is
isomorphic to $(\bE^\wedge_\cJ,\cT^\wedge_\cJ)$, and the
completion $\bE_{I+J}^\wedge$ of $\bE(A,\cT_{I+J})$ is
isomorphic to $(\bE^\wedge_{\cI+\cJ},\cT^\wedge_{\cI+\cJ})$.

Next, set $\cK:=\cI\cJ+\alpha_0\bE\subset\bE$; since $\cI$
and $\cJ$ are finitely generated and they both contain
$\alpha_0$, it is easily seen that the $(\cI\cap\cJ)$-adic
topology on $\bE$ agrees with the $\cK$-adic topology.
Likewise, topology $\cT_{I\cap J}$ on $A$ agrees with the
$(IJ)$-adic topology. Moreover, $\cK$ is a finitely
generated ideal, by construction $p\in IJ$, and
$\bar u_{A/pA}(\cK)=(IJ)/pA$, hence remark
\ref{rem_nice-topology}(iii) also tells us that
$(\bE^\wedge_{\cI\cap\cJ},\cT^\wedge_{\cI\cap\cJ})$ is isomorphic
to the completion $\bE_{I\cap J}^\wedge$ of $\bE(A,\cT_{I\cap J})$.

Notice now that, since the image of $A$ is dense in
$A^\wedge_{I\cap J}$, the map $\phi_I$ in
\eqref{eq_old-acquaintance} is the unique continuous
homomorphism of topological $A$-algebras from
$A^\wedge_{I\cap J}$ to $A^\wedge_I$, and the other
maps in \eqref{eq_old-acquaintance} enjoy corresponding
uniqueness properties. Combining with remark
\ref{rem_nice-topology}(iv), we conclude that
$W(\cE)\otimes_\bA A$ is naturally identified with
\eqref{eq_old-acquaintance}, and the image of
$\underline\alpha$ is still distinguished in
each of the Witt rings appearing in $W(\cE)$.
Due to of claim \ref{cl_still-cart}, we already
see that \eqref{eq_old-acquaintance} is a
cocartesian diagram of $A$-modules;
moreover, it will follow that it is a cartesian
diagram of rings, once we know that
\set\begin{equation}\label{eq_standard-calc}
\Tor^\bA_1(W(\bE^\wedge_{\cI+\cJ}),A)=0.
\end{equation}
However $A=\bA/\underline\alpha\bA$, and
$\underline\alpha$ is a regular element of
$W(\bE^\wedge_{\cI+\cJ})$ (remark \ref{rem_distinguished}(ii)),
so \eqref{eq_standard-calc} holds by a standard calculation.

(ii): We know already that \eqref{eq_double-completion}
is open and surjective (lemma \ref{lem_double-completion}(i)),
so it suffices to show that this map is a ring isomorphism.
However, in light of remark \ref{rem_nice-topology}(iv)
we have a diagram of continuous maps
$$
\cD\quad :\quad
{\diagram (A^\wedge_I)^\wedge_J \ar[r] \ar[d] &
A(\bE^\wedge_\cI,\alpha_0\bE^\wedge_\cI)^\wedge_J \ar[r] &
A((\bE^\wedge_\cI)^\wedge_\cJ,\alpha_0(\bE^\wedge_\cI)^\wedge_\cJ)
\ar[d] \\
A^\wedge_{I+J} \ar[rr] & &
A(\bE^\wedge_{\cI+\cJ},\alpha_0\bE^\wedge_{\cI+\cJ})
\enddiagram}$$
whose left vertical arrow equals \eqref{eq_double-completion},
and whose remaining arrows are ring isomorphisms. Thus,
it suffices to check that $\cD$ commutes, {\em i.e.}
that the two continuous ring homomorphisms
$f,g:(A^\wedge_I)^\wedge_J\to
A(\bE^\wedge_{\cI+\cJ},\alpha_0\bE^\wedge_{\cI+\cJ})$ deduced
from $\cD$ coincide. However, let $i:A\to (A^\wedge_I)^\wedge_J$
denote the completion map; a simple inspection shows that
$f\circ i=g\circ i$, and since the topology of
$A(\bE^\wedge_{\cI+\cJ},\alpha_0\bE^\wedge_{\cI+\cJ})$ is
separated and $i$ has dense image, the assertion follows. 
\end{proof}

\sset\subsubsection{}\label{subsec_construct-new-perfs}
Let $A$, $\bar D$, $\bar A$ be three perfectoid rings, and set
$\bE:=\bE(A)$, $\bar\bE:=\bE(\bar D)$. Suppose that
$\bar A$ is a discrete topological $\F_p$-algebra, so
that $\bE(\bar A)=\bar A$ (see also corollary
\ref{cor_perf-are-reduced}(iii)), let
$\bar\phi_A:\bar D\to\bar A$ and $\pi_A:A\to\bar A$ be two
continuous ring homomorphisms, and suppose as well that
$\pi_A$ is surjective. Define the topological rings $D_A$
and $D_\bE$ as the fibre products in the resulting cartesian
diagrams of topological rings
$$
\xymatrix{ D_A \ar[r]^-{\phi_A} \ar[d] & A \ar[d]^{\pi_A} &
D_\bE \ar[r]^-{\phi_\bE} \ar[d] &
\bE \ar[d]^{\pi_\bE} \\
\bar D \ar[r]^-{\bar\phi_A} & \bar A &
\bar\bE \ar[r]^-{\bar\phi_\bE} & \bar A.
}$$
where $\pi_\bE:=\bE(\pi_A)$ and $\bar\phi_\bE:=\bE(\bar\phi_A)$.

\begin{proposition}\label{prop_construct-new-perfs}
In the situation of \eqref{subsec_construct-new-perfs},
the rings $D_A$ and $D_\bE$ are perfectoid, and there
exists a natural isomorphism of topological rings
$$
\omega:\bE(D_A)\isom D_\bE
\qquad\text{such that}\qquad
\phi_\bE\circ\omega=\bE(\phi_A).
$$
\end{proposition}
\begin{proof} To begin with, we notice that $D_A$ (resp.
$D_\bE$) is complete and separated, since $\bar A$ is
separated and both $\bar D$ and $A$ (resp. $\bar\bE$
and $\bE$) are complete and separated. Let
$(\alpha_n~|~n\in\N)$ (resp. $(\alpha'_n~|~n\in\N)$)
be a distinguished element in $\Ker\,u_A$ (resp. in
$\Ker\,u_{\bar D}$), and pick a finitely generated ideal
of adic definition $J_\bE\subset\Ker\,\pi_\bE$
(resp. $J_{\bar\bE}\subset\Ker\,\bar\phi_\bE$)
for $\bE$ (resp. for $\bar\bE$).

\begin{claim}\label{cl_exp-3}
We may assume that $\alpha^{1/p^2}_0\bE\subset J_\bE$ and
$\alpha'^{1/p^2}_0\bar\bE\subset J_{\bar\bE}$.
\end{claim}
\begin{pfclaim} For the first stated inclusion, notice
that $\alpha_0^{1/p^n}\in\Ker\,\pi_\bE$ for every $n\in\N$,
so $J_\bE+\alpha_0^{1/p^n}\bE$ is still an ideal of adic
definition of $\bE$ contained in $\Ker\,\pi_\bE$, for
every such $n$. The same argument applies to $ J_{\bar\bE}$.
\end{pfclaim}

By virtue of claim \ref{cl_exp-3} we may find a system of
generators $\beta_\bullet:=(\beta_1,\dots,\beta_k)$ for $J_\bE$
with $\beta_1=\alpha_0^{1/p^2}$, and we set
$b_i:=\bar u_A(\beta_i)$ for $i=1,\dots,k$. It
follows that $b_\bullet:=(b_1,\dots,b_k)$ is a system of
generators for an ideal $J_A$ of adic definition of $A$.
Likewise, pick a system of generators
$\beta'_\bullet:=(\beta'_1,\dots,\beta'_h)$ for $J_{\bar\bE}$
with $\beta'_1=\alpha'^{1/p^2}_0$, and set
$b'_i:=\bar u_{\bar D}(\beta'_1)$ for $i=1,\dots,h$. Then
$b'_\bullet:=(b'_1,\dots,b'_h)$ is a system of generators
for an ideal $J_{\bar D}$ of adic definition of $\bar D$.
Notice that $J_A\times J_{\bar D}\subset D_A$ (resp.
$J_\bE\times J_{\bar\bE}\subset D_\bE$), and the topology
of $D_A$ (resp. of $D_\bE$) is the linear topology defined
by the system of ideals $(J^n_A\times J^n_{\bar D}~|~n\in\N)$
(resp. $(J^n_\bE\times J^n_{\bar\bE}~|~n\in\N)$).
Denote by $K_\bE\subset D_\bE$ the ideal generated
by $(\beta_\bullet\times\{0\},\{0\}\times\beta'_\bullet)$;
it is easily seen that
$$
J^{n+1}_\bE\times J^{n+1}_{\bar\bE}\subset
K^n_\bE\subset J^n_\bE\times J^n_{\bar\bE}
\qquad
\text{for every $n\in\N$}
$$
so the topology of $D_\bE$ agrees with its $K_\bE$-adic
topology, and therefore $D_\bE$ is a perfectoid $\F_p$-algebra.
Likewise, let $K_A\subset D_A$ be the ideal generated
by $(b_\bullet\times\{0\},\{0\}\times b'_\bullet)$; by
the same token, we have
$$
J^{n+1}_A\times J^{n+1}_{\bar D}\subset
K^n_A\subset J^n_A\times J^n_{\bar D}
\qquad
\text{for every $n\in\N$}
$$
so the topology of $D_A$ agrees with the $K_A$-adic topology.
Moreover, by construction we have $pA\subset J_A^{p^2}$ and
$p\bar D\subset J_{\bE}^{p^2}$ (lemma \ref{lem_was-third-cond}(iii));
since $p^2-1\geq 2$, we deduce that $pD_A\subset K_A^2$, and
we see already that $D_A$ is a P-ring. Next, by remark
\ref{rem_Witt-limit}(ii), we deduce a short exact sequence
of $W(D_A)$-modules
$$
\cW
\quad :\quad
0\to W(D_A)\to W(A)\oplus W(\bar D)\to W(\bar A)\to 0.
$$
On the other hand, by lemma \ref{lem_was-third-cond}(ii)
we may find a distinguished element
$\underline\alpha'':=(\alpha''_n~|~n\in\N)$ in
$\Ker\,u_{D_A}$ and to ease notation we set
$D'_A:=W(D_A)/\underline\alpha'' W(D_A)$; since the
image of $\underline\alpha''$ is still distinguished in
$W(\bar A)$, remark \ref{rem_distinguished}(ii) implies
that $\Tor^{W(D_A)}_1(D'_A,W(\bar A))=0$, so the sequence
$$
\cW\otimes_{W(D_A)}D'_A
\quad :\quad
0\to D'_A\to A\oplus\bar D\to\bar A\to 0
$$
is still exact, and consequently $u_{D_A}$ induces a ring
isomorphism $D'_A\to D_A$, {\em i.e.} $D_A$ is perfectoid.
Lastly, notice that $u_A\otimes_\Z\F_p$ and $u_{\bar D}$
induce isomorphisms of topological $\F_p$-algebras
$$
D_A/pD_A\isom\bar D/p\bar D\times_{\bar A}A/pA\isom
\bar D/\alpha''_0\bar D\times_{\bar\bE}\bE/\alpha''_0\bE
\isom D_\bE/\alpha''_0D_\bE.
$$
Combining with theorem \ref{th_fontaine} we deduce an
isomorphism
$$
\bE(D_A)\isom\bE(D_A/pD_A)\isom\bE(D_\bE/\alpha''_0D_\bE)\isom D_\bE
$$
which fulfills the stated condition, by a simple inspection.
\end{proof}

\begin{definition}\label{def_beta-taut}
Let $A$ be a perfectoid ring and
$\cJ\subset\cK\subset\bE:=\bE(A)$ any two ideals.
Let also $\underline\alpha:=(\alpha_n~|~n\in\N)\in\Ker\,u_A$
be any distinguished element and $\beta\in\bE$ any element.
\begin{enumerate}
\item
We say that the inclusion of $\cJ$ in $\cK$ is
{\em $\beta$-taut} if we have
$$
\beta\cdot\Phi^{-1}_\bE(\cK^p)\subset\cJ.
$$
\item
We say that $\cJ$ is {\em $\beta$-taut} if the
identity map of $\cJ$ is a $\beta$-taut inclusion.
\item
We say that $\cJ$ is {\em strictly $\beta$-taut}
if it is $\beta^\lambda$-taut for some
$\lambda\in\Z[1/p]$ with $0\leq\lambda<1$.
\item
$\cJ$ is {\em taut} (resp. {\em strictly taut}) if
it is $\alpha_0$-taut (resp. strictly $\alpha_0$-taut),
and the inclusion $\cJ\subset\cK$ is {\em taut} (resp.
{\em strictly taut}), if it is $\alpha_0$-taut (resp.
strictly $\alpha_0$-taut).
\item
We denote by $\{\cJ\}\subset A$ the topological closure of
the ideal generated by the system $(\bar u_A(x)~|~x\in\cJ)$
(notation of \eqref{subsec_Teich}).
\end{enumerate}
\end{definition}

\begin{remark}\label{rem_beta-taut}
With the notation of definition \ref{def_beta-taut},
the following holds.

(i)\ \
It follows easily from remark \ref{rem_distinguished}(i)
that the definition of taut and strictly taut ideals does
not depend on the choice of $\underline\alpha$.

(ii)\ \
By inspecting the definition of angular powers, it is
easily seen that $\cJ$ is $1$-taut if and only if
$\cJ=\cJ^{\La 1\Ra}\bE$.

(iii)\ \
It follows from (ii) and lemma \ref{lem_mon-fract-powers}(ii.c)
that for every ideal $\cJ$ of $\bE$ and every 
$\lambda\in\N[1/p]$, the ideal $\cJ^{\La\lambda\Ra}\bE$
is strictly taut. Moreover, say that
$\alpha_0\in\cJ^{\La\eps\Ra}\bE$ for some $\eps\in\N[1/p]$;
then the inclusion
$\cJ^{\La\lambda+\eps\Ra}\bE\subset\cJ^{\La\lambda\Ra}\bE$ is
taut for every $\lambda\in\N[1/p]$, by lemma
\ref{lem_mon-fract-powers}(ii.b).

(iv)\ \
If $\cJ_1,\cJ_2\subset\bE$ are two $\beta$-taut
(resp. strictly $\beta$-taut, resp. taut, resp.
strictly taut) ideals, then the same holds for
$\cJ_1\cap\cJ_2$. Indeed, suppose that both
ideals are $\beta$-taut; then we have :
$$
\begin{aligned}
\beta\cdot\Phi_\bE^{-1}((\cJ_1\cap\cJ_2)^p)\subset\, &
\beta\cdot\Phi_\bE^{-1}(\cJ_1^p\cap\cJ_2^p) \\
=\, &
\beta\cdot(\Phi_\bE^{-1}(\cJ_1^p)\cap\Phi_\bE^{-1}(\cJ_2^p)) \\
\subset\, & (\beta\cdot\Phi_\bE^{-1}(\cJ_1^p))\cap
(\beta\cdot\Phi_\bE^{-1}(\cJ_2^p)) \\
\subset\, & \cJ_1\cap\cJ_2
\end{aligned}
$$
whence the claim. In the same vein, if $\cJ_1\subset\cJ_2$
is a $\beta$-taut inclusion, and $\cK$ is any $\beta$-taut
ideal of $\bE$, then the inclusion
$\cK\cap\cJ_1\subset\cK\cap\cJ_2$ is $\beta$-taut as well
(details left to the reader); the same holds if $\beta$-taut
is replaced by strictly $\beta$-taut, taut or strictly taut.

(v)\ \
If $\cJ$ is $\beta$-taut, then the same holds
for the topological closure $\cJ^c$ of $\cJ$ in $\bE$.
Likewise, if $\cJ\subset\cK$ is a $\beta$-taut
inclusion, then the same holds for the inclusion
$\cJ^c\subset\cK^c$ of the respective topological
closures; moreover, both $\cJ$ and $\cK$ are
$\beta$-taut. Then clearly the same assertions hold
with $\beta$-taut replaced by strictly $\beta$-taut,
taut, or strictly taut (details left to the reader).
In the same vein, notice that
$$
\{\cJ^c\}=\{\cJ\}
\qquad
\text{for every ideal $\cJ\subset\bE$}.
$$

(vi)\ \
Let $\cJ_1,\cJ_2,\cJ_3\subset\bE$ be three ideals,
$n,m\in\N$ any two integers, and suppose that
$$
\beta_1\cdot\Phi^{-n}_\bE(\cJ_1^{p^n})\subset\cJ_2
\quad\text{and}\quad
\beta_2\cdot\Phi^{-m}_\bE(\cJ_2^{p^m})\subset\cJ_3
\qquad
\text{for some $\beta_1,\beta_2\in\bE$}.
$$
Then it is easily seen that
$\beta_1\beta_2\cdot\Phi^{-n-m}_\bE(\cJ_1^{p^{n+m}})\subset\cJ_3$.
Especially, if $\cJ$ is $\beta$-taut, we have
\set\begin{equation}\label{eq_iterate-beta}
\beta^n\cdot\Phi^{-n}_\bE(\cJ^{p^n})\subset\cJ
\qquad
\text{for every $n\in\N$}.
\end{equation}

(vii)\ \
For $i=1,2$, let $\beta_i\in\bE$ be any element, and
$\cJ_i$ a $\beta_i$-taut ideal; then it is easily seen
that $\cJ_1\cJ_2$ is $\beta_1\beta_2$-taut. More generally,
if $\cJ\subset\cK$ is a $\beta_1$-taut inclusion of ideals
of $\bE$, and $\cJ'\subset\bE$ is any other $\beta_2$-taut
ideal, then the inclusion $\cJ\cJ'\subset\cK\cJ'$ is
$\beta_1\beta_2$-taut.

(viii)\ \
Suppose that $\cJ$ is an open ideal of $\bE$.
Then the system $(\bar u_A(x)~|~x\in\cJ)$ generates
an open ideal $J$ of $A$, and therefore $\{\cJ\}=J$ is
open as well in $A$. Indeed, let $I\subset A$ be any
ideal of definition, pick a finite system
$(\bar a_1,\dots,\bar a_k)$ of generators of $I/pA$,
let $\alpha_1,\dots,\alpha_k$ be elements
of $\bE$ such that $\bar u_{A/pA}(\alpha_i)=\bar a_i$
for $i=1,\dots,k$, and denote by $\cI\subset\bE$ the
ideal generated by the system
$(\alpha_1,\dots,\alpha_k)$. Arguing as in the proof
of lemma \ref{lem_was-third-cond}(i), we see that the
system $(\bar u_A(\alpha_i)~|~i=1,\dots,k)$ generates
$I$, and $\cI^n\subset\cJ$ for every sufficiently large
$n\in\N$, so that $I^n\subset J$ for every such $n$,
whence the assertion.
\end{remark}

\sset\subsubsection{}\label{subsec_elementary}
Keep the notation of definition \ref{def_beta-taut}, and
let $\cJ\subset\cK$ be a taut inclusion of ideals of
$\bE$. Notice that $\alpha_0\cK\subset\cJ$, so $\cK/\cJ$
is an $\bE/\alpha_0\bE$-module, and then it can be viewed
as an $A/pA$-module, via the isomorphism $\omega$ of
remark \ref{rem_nice-topology}(ii). Likewise, from lemma
\ref{lem_was-third-cond}(iii) we easily deduce that
$p\{\cK\}\subset\{\cJ\}$, so $\{\cK\}/\{\cJ\}$ is an
$A/pA$-module as well.

\begin{theorem}\label{th_taut-one}
Let $A$ be any perfectoid ring, and set\/
$\bE:=\bE(A)$. The following holds :
\begin{enumerate}
\item
Every taut inclusion $\cJ_1\subset\cJ_2$ of ideals of\/
$\bE$ induces an $A/pA$-linear map
$$
\bar\tau:\cJ_2/\cJ_1\to\{\cJ_2\}/\{\cJ_1\}
\quad :\quad
(x\mod{\cJ_1})\mapsto(\bar u_A(x)\mod{\{\cJ_1\}}).
$$
\item
If both $\cJ_1$ and $\cJ_2$ are closed in the topology
of\/ $\bE$, the map\/ $\bar\tau$ of\/ {\em (i)} is an
isomorphism.
\end{enumerate}
\end{theorem}
\begin{proof} Let
$\underline\alpha:=(\alpha_n~|~n\in\N)\in\Ker\,u_A$ be
any distinguished element.

(i): Let $x,y\in\cJ_2$ be any two elements, and set
$\delta:=\tau_A(x+y)-\tau_A(x)-\tau_A(y)$. In light of
lemma \ref{lem_drop-conditions}(iv), we have to check
that $u_A(\delta)$ lies in $\{\cJ_1\}$.
However, proposition \ref{prop_combinatorial} expresses
$\delta$ as a series of the form $\sum_{n\in\N}p^nb_n$,
where each $b_n$ is in turn a finite sum of terms of the
form $c_{n,\sigma}\tau_A(\beta_{n,\sigma})$ (for $\sigma$
ranging over a certain finite set $\Sigma_n$), and
with $\beta_{n,\sigma}\in\Phi_\bE^{-n}(\cJ_2^{p^n})$ and
$c_{n,\sigma}\in\Z_p$ for every $n\in\N$ and every
$\sigma\in\Sigma_n$. Thus, we come down to checking
that $p^n\cdot\bar u_A(\beta_{n,\sigma})\in\{\cJ_1\}$
for every $n\in\N$ and $\sigma\in\Sigma_n$. In light of
lemma \ref{lem_was-third-cond}(iii), this holds if and only
if $\bar u_A(\alpha_0^n\cdot\beta_{n,\sigma})\in\{\cJ_1\}$,
which in turns will follow, once we know that
$\alpha_0^n\cdot\beta_{n,\sigma}\in\cJ_1$. But the
latter is clear from remark \ref{rem_beta-taut}(vi)
(with $\beta:=\alpha_0$ : details left to the reader).

Lastly, for any $x\in\cJ_2$ and $\beta\in\bE$ we have
$\tau_A(\beta x)=\tau_A(\beta)\cdot\tau_A(x)$, and the
isomorphism $\omega$ of remark \ref{rem_nice-topology}(ii)
maps the class of $\beta$ in $\bE/\alpha_0\bE$ to the class
of $\bar u_A(\beta)$ in $A/pA$; this shows that
$\bar\tau(\bar a\cdot x)=\bar a\cdot\bar\tau(x)$ for
every $\bar a\in A/pA$ and every $x\in\cJ_2/\cJ_1$,
whence (i).

Next, choose a strictly positive $\eps\in\N[1/p]$
such that $\alpha_0\in\cJ^{\La\eps\Ra}\bE$; we point out
the following special case of (ii) :

\begin{claim}\label{cl_special-case}
Let $\cJ$ be any ideal of definition of
$\bE$, and $r,r'\in\R_+$ two real numbers such that
$\eps\geq r'-r\geq 0$. Then (ii) holds for the taut
inclusion
$\cJ^{\lfloor r'\rfloor}\bE\subset\cJ^{\lfloor r\rfloor}\bE$.
\end{claim}
\begin{pfclaim} Let us set
$W\lfloor s\rfloor:=W(\cJ^{\lfloor s\rfloor}\bE)$
for every $s\in\N[1/p]$ (notation of remark
\ref{rem_Witt-limit}(iv)); from \eqref{eq_ghost-angular}
we see that $\tau_A$ induces a natural isomorphism
$$
\cJ^{\lfloor r\rfloor}\bE/\cJ^{\lfloor r'\rfloor}\bE\isom
W\lfloor r\rfloor/(pW\lfloor r\rfloor+W\lfloor r'\rfloor)
$$
and we are reduced to showing that the restriction
of $u_A$ to $W\lfloor r\rfloor$ induces an isomorphism
\set\begin{equation}\label{eq_W-and-brackets}
W\lfloor r\rfloor/(pW\lfloor r\rfloor+W\lfloor r'\rfloor)
\isom\{\cJ^{\lfloor r\rfloor}\bE\}/\{\cJ^{\lfloor r'\rfloor}\bE\}.
\end{equation}
However, remark \ref{rem_beta-taut}(viii) already
implies that \eqref{eq_W-and-brackets} is surjective,
and its kernel is the image of
$W\lfloor r\rfloor\cap\underline\alpha W(E)$. In light of
proposition \ref{prop_square-powers}(ii) we then come down
to showing :
$$
\underline\alpha W\lfloor r\rfloor\subset
pW\lfloor r\rfloor+W\lfloor r'\rfloor.
$$
But due to our choice of $\eps$ we get
$\alpha_0\cJ^{\lfloor r\rfloor}\bE\subset\cJ^{\lfloor r'\rfloor}\bE$,
so the latter inclusion follows again from
\eqref{eq_ghost-angular}.
\end{pfclaim}

\begin{claim}\label{cl_little-top-trick}
Let $B$ be any ring, $f:M\to M'$ and $g:M'\to M''$ two
continuous maps of topological $B$-modules, and suppose
that
\begin{enumerate}
\alphaenu
\item
the topologies of $M$ and $M''$ are discrete and
that of $M'$ is separated
\item
$f$ has dense image and $g\circ f$ is injective.
\end{enumerate}
Then $f$ is an isomorphism of topological $B$-modules.
\end{claim}
\begin{pfclaim} By assumption, the topological
closure of $f(M)$ in $M'$ equals $M'$; by claim
\ref{cl_image-and-closure} it follows that the
topological closure of $g\circ f(M)$ in $M''$
equals the topological closure of $g(M')$.
But since the topology of $M''$ is discrete, this
just means that $g\circ f(M)=g(M')$. We may therefore
replace $M''$ by $g(M')$, and assume from start
that $g$ is surjective and $g\circ f$ is an
isomorphism. Then we may even assume that $M''=M$
and $g\circ f=\one_M$. Let us endow $\Img\,f$ and
$\Ker\,g$ with the topologies induced from the
inclusion into $M'$; then the addition law of $M'$
restricts to a continuous and bijective $B$-linear map
$$
h:\Img\,f\oplus\Ker\,g\to M'
$$
(where the direct sum is endowed with the product
topology). However, the inverse of $h$ is the map
given by the rule :
$m'\mapsto(f\circ g(m'),m'-f\circ g(m'))$ for every
$m'\in M'$. Clearly, this map is also continuous,
so $h$ is an isomorphism of topological $B$-modules.
Now, since $M'$ is separated, the same must hold
for $\Ker\,g$; especially, $L:=h(\Img\,f\oplus 0)$ is
a closed subset of $M'$. On the other hand,
$\Img\,f=L$; since $f$ has dense image, we must then
have $L=M'$, {\em i.e.} $\Ker\,g=0$, whence the claim.
\end{pfclaim}

We may now complete the proof of (ii): fix an ideal
of definition $\cJ$ of $\bE$ with $\alpha_0\in\cJ$,
and for every ideal $\cK\subset\bE$ set
$$
\Fil^n\cK:=\cK\cap\cJ^{\lfloor n\rfloor}\bE
\qquad
\Fil^n\{\cK\}:=\{\Fil^n\cK\}
\qquad
\text{for every $n\in\N$}.
$$
If $\cK'\subset\cK$ is any inclusion of ideals of $\bE$,
the image of $\Fil^\bullet\cK$ in $\cK/\cK'$ defines a
filtration $\Fil^\bullet(\cK/\cK')$ of $\cK/\cK'$, and
likewise we get a filtration $\Fil^\bullet(\{\cK\}/\{\cK'\})$
on $\{\cK\}/\{\cK'\}$. Notice that the inclusion
$\Fil^{n+1}\cK\subset\Fil^n\cK$ is taut for every $n\in\N$
and every taut ideal $\cK$, and denote by $\gr^\bullet\cK$
(resp. $\gr^\bullet\{\cK\}$) the graded $A/pA$-module
associated to the filtration $\Fil^\bullet\cK$ (resp.
to $\Fil^\bullet\{\cK\}$). Likewise we define
$\gr^\bullet(\cK/\cK')$ and $\gr^\bullet(\{\cK\}/\{\cK'\})$
for an inclusion of ideals $\cK'\subset\cK$. We remark :

\begin{claim}\label{cl_induction-lambda}
If $\cK$ is taut, the map $\bar\tau$ of (i) induces an
isomorphism
$$
\gr^n\cK\isom\gr^n\{\cK\}
\qquad
\text{for every $n\in\N$}.
$$
\end{claim}
\begin{pfclaim} Let us endow any ideal of $\bE$ (resp. of $A$)
with the topology induced from the inclusion into $\bE$ (resp.
into $A$), and for any inclusion $\cI\subset\cI'$ (resp.
$I\subset I'$) of ideals of $\bE$ (resp. of $A$), let us endow
$\cI'/\cI$ (resp. $I'/I$) with the corresponding quotient topology.
Notice that the inclusions $\Fil^{n+1}\cK\subset\Fil^n\cK$
and $\Fil^{n+1}\{\cK\}\subset\Fil^n\{\cK\}$ are both taut
(remark \ref{rem_beta-taut}(iv)); from (i) we deduce a
commutative diagram of $A/pA$-modules :
$$
\xymatrix{
\gr^n\cK \ar[r] \ar[d] & \gr^n\{\cK\}
\ar[d] \\
\gr^n\bE \ar[r] & \gr^nA
}$$
whose left vertical arrow is injective, and whose bottom
horizontal arrow is already known to be an isomorphism,
by claim \ref{cl_special-case}.
Moreover, all these maps are continuous, for the topologies
that we have just defined on these modules, and furthermore,
the top horizontal arrow has dense image. Also, since
$\cJ$ is open, the same holds for $\cJ^{\lfloor r'\rfloor}\bE$,
and therefore the two modules on the bottom row have
both the discrete topology. Likewise,
$\cK\cap\cJ^{\lfloor r'\rfloor}\bE$ is open in
$\cK\cap\cJ^{\lfloor r'\rfloor}\bE$, so also the source of
the map on the top row is a discrete $A/pA$-module.
Lastly, since $\{\cK\cap\cJ^{\lfloor r'\rfloor}\bE\}$ is
a closed ideal, the target of the same map is a separated
$A/pA$-module. Thus, all the conditions of claim
\ref{cl_little-top-trick} are fulfilled, and it follows
that the top horizontal arrow is an isomorphism.
\end{pfclaim}

There follows, for every $n\in\N$ a commutative diagram
with exact rows
$$
\xymatrix{
0 \ar[r] & \gr^n\cJ_1 \ar[r] \ar[d] &
\gr^n\cJ_2 \ar[r] \ar[d] &
\gr^n(\cJ_2/\cJ_1) \ar[r] \ar[d] & 0 \\
0 \ar[r] & \gr^n\{\cJ_1\} \ar[r] & \gr^n\{\cJ_2\} \ar[r] &
\gr^n(\{\cJ_2\}/\{\cJ_1\}) \ar[r] & 0
}$$
and claim \ref{cl_induction-lambda} shows that the leftmost
and central vertical arrows are both isomorphisms, so the
same holds for the rightmost vertical arrow. Now, notice
that the filtration $\Fil^\bullet(\cJ_2/\cJ_1)$ defines a
separated and complete topology on $\cJ_2/\cJ_1$, since
$\cJ_2$ and $\cJ_1$ are both closed ideals in $\bE$.
The same holds for the topology on $\{\cJ_2\}/\{\cJ_1\}$
determined by the filtration
$\Fil^\bullet(\{\cJ_2\}/\{\cJ_1\})$, since $\{\cJ_2\}$
and $\{\cJ_1\}$ are closed ideals in $A$. Then (ii)
follows directly from \cite[Ch.III, \S2, n.8, Cor.3]{BouAC}.
\end{proof}

\begin{theorem}\label{th_taut-two}
Let $A$ be any perfectoid ring, $\cK$ a taut ideal of\/
$\bE:=\bE(A)$, and $W(\cK)^c$ the topological closure
of\/ $W(\cK)$ in $\bA(A)$ (notation of remark
{\em\ref{rem_Witt-limit}(iv)}). We have :
\begin{enumerate}
\item
$\{\cK\}=u_A(W(\cK)^c)$.
\item
$\cK$ is an open ideal if and only if the same holds
for $\{\cK\}$.
\item
If $\cK'\subset\cK$ is another taut ideal, and both
$\cK$ and $\cK'$ are closed, then $\cK=\cK'$ if and
only if $\{\cK\}=\{\cK'\}$.
\item
Suppose that $\cK$ is $1$-taut and closed, and let
$\underline\alpha$ be any distinguished element of\/
$\Ker\,u_A$. Then the pair $(\underline\alpha,\cK)$
is transversal (see \eqref{subsec_transvesal}).
\end{enumerate}
\end{theorem}
\begin{proof} Pick any distinguished element
$\underline\alpha:=(\alpha_n~|~n\in\N)$ of $\Ker\,u_A$.

(i): First we prove that $u_A(W(\cK)^c)\subset\{\cK\}$.
Indeed, let $\underline a:=(a_n~|~n\in\N)\in W(\cK)$ be
any element; taking into account \eqref{eq_new-form},
we are easily reduced to showing that
$$
u_A(p^n\cdot\tau_A(a^{p^{-n}}_n))\in\{\cK\}
\qquad
\text{for every $n\in\N$}.
$$
In light of lemmata \ref{lem_was-third-cond}(iii) and
\ref{lem_drop-conditions}(iv), it then suffices to check
that $\bar u_A(\alpha_0^n\cdot a^{p^{-n}}_n)\in\{\cK\}$,
and the latter is clear from \eqref{eq_iterate-beta}
(applied with $\beta:=\alpha_0$). Next, fix an ideal
of definition $\cJ$ of $\bE$ containing $\alpha_0$,
and define the filtrations $\Fil^\bullet\cK$ on $\cK$
and $\Fil^\bullet\{\cK\}$ on $\{\cK\}$ as in the proof
of theorem \ref{th_taut-one}(ii); by
\eqref{eq_general-inclusion} and remark
\ref{rem_beta-taut}(v) we may assume that $\cK$ is
a closed ideal of $\bE$, in which case theorem
\ref{th_taut-one}(ii) shows that the map $\bar u_A$
induces an isomorphism $\gr^\bullet\cK\isom\gr^\bullet\{\cK\}$
on the respective associated graded $A/pA$-modules.
Now, let $x\in\{\cK\}$ be any element; it follows
easily that we may write
$$
x=u_A\Bigl(\sum_{n\in\N}\tau_A(\beta_n)\Bigr)
\qquad
\text{where $\beta_n\in\Fil^n\cK$ for every $n\in\N$}
$$
and the series converges in the topology of $\bA(A)$.
But clearly $\tau_A(\beta_n)\in W(\cK)$ for every
$n\in\N$, whence the sought converse inclusion.

(ii): In light or remark \ref{rem_beta-taut}(viii) we may
assume that $\{\cK\}$ is open, and we need to show
that the same holds for $\cK$; then, by lemma
\ref{lem_5.3.8}(ii.b), we may also assume that
$\cK$ is a closed ideal. The inclusion of $\cK$ into
$\bE$ is a map of filtered $\bE$-modules
$\Fil^\bullet\cK\to\Fil^\bullet\bE$, and likewise
we get a map of filtered $A$-modules
$\Fil^\bullet\{\cK\}\to\Fil^\bullet A$. The assumption
on $\{\cK\}$ implies that there exists $n\in\N$ such
that the map of associated graded modules
$\gr^k\{\cK\}\to\gr^kA$ is an isomorphism for every
$k\geq n$. Then theorem \ref{th_taut-one}(ii) implies
that the same holds for the corresponding map
$\gr^k\cK\to\gr^k\bE$. However, both $\cK$ and $\bE$
are complete and separated for their filtrations, so
$\Fil^n\cK=\Fil^n\bE=\cJ^{\La n\Ra}$
(\cite[Ch.III, \S2, n.8, Cor.3]{BouAC}), {\em i.e.}
$\cK$ is open.

(iii) is similar : we may assume that $\{\cK\}=\{\cK'\}$,
and we consider the induced map of filtered $\bE$-modules
$\Fil^\bullet\cK'\to\Fil^\bullet\cK$. Arguing as in the
proof of (ii), we see that the latter induces an
isomorphism on the associated graded $\bE$-modules, so
$\cK=\cK'$, again by \cite[Ch.III, \S2, n.8, Cor.3]{BouAC}.

(iv): Denote by  $(\alpha_0\cK)^c$ the topological closure
of $\underline\alpha\cK$ in $\bE$. We shall show first that

\begin{claim}\label{cl_start-here}
$W(\cK)\cap\Ker\,u_A=W((\alpha_0\cK)^c)\cap\Ker\,u_A
+\underline\alpha W(\cK)$.
\end{claim}
\begin{pfclaim} Indeed, let
$\underline\omega:=(\omega_n~|~n\in\N)\in W(\cK)\cap\Ker\,u_A$
be any element; due to remarks \ref{rem_Witt-are-f-adic}(iv)
and \ref{rem_beta-taut}(ii), we may write
$$
\underline\omega=\tau_A(\omega_0)+p\cdot\underline\omega'
\qquad
\text{for some $\underline\omega'\in W(\cK)$}
$$
whence $0=u_A(\underline\omega)=\bar u_A(\omega_0)+
p\cdot u_A(\underline\omega')$, and especially, the image
of $\bar u_A(\omega_0)$ vanishes in $\{\cK\}/\{(\alpha_0\cK)^c\}$
(lemma \ref{lem_was-third-cond}(iii)). In view of theorem
\ref{th_taut-one}(ii), we deduce that
$\omega_0\in(\alpha_0\cK)^c$. Write also
$\underline\alpha=\tau_A(\alpha_0)+p\underline u$ for some
$\underline u\in\bA(A)^\times$; we obtain
$$
\underline\omega=\tau_A(\omega_0)+
(\underline u^{-1}\underline\alpha-\tau_A(\alpha_0))
\cdot\underline\omega')
$$
and it suffices to notice that
$\underline u^{-1}\underline\alpha\cdot\underline\omega'\in
\underline\alpha W(\cK)\subset W(\cK)\cap\Ker\,u_A$, and
consequently
$\tau_A(\omega_0)-\tau_A(\alpha_0)\cdot\underline\omega'\in
W((\alpha_0\cK)^c)\cap\Ker\,u_A$.
\end{pfclaim}

Notice next that if $\cK$ is $1$-taut, then the same holds
for $\alpha_0\cK$, and then also for $(\alpha_0\cK)^c$
(remark \ref{rem_beta-taut}(v)). Thus, we may apply claim
\ref{cl_start-here} to $(\alpha^n_0\cK)^c$, and get
$$
W((\alpha^n_0\cK)^c)\cap\Ker\,u_A=W((\alpha^{n+1}_0\cK)^c)\cap\Ker\,u_A
+\underline\alpha W((\alpha^n_0\cK)^c)
\qquad
\text{for every $n\in\N$}.
$$
Now, let $x_0\in W(\cK)\cap\underline\alpha W(\bE)$ be
any element; we may then find inductively identities
$$
x_n=\underline\alpha\cdot y_n+x_{n+1}
\qquad
\text{with $x_{n+1}\in W((\alpha^{n+1}\cK)^c)\cap\Ker\,u_A$
and $y_n\in W((\alpha_0^n\cK)^c)$}
$$
for every $n\in\N$. Therefore :
$$
x_0=\underline\alpha\cdot(y_0+y_1+\cdots+y_n)+x_{n+1}
\qquad
\text{for every $n\in\N$}.
$$
Lastly, since $\cK$ is $1$-taut, the ideal $W(\cK)$ is
closed in the topology of $W(\bE)$ (see remark
\ref{rem_Witt-limit}(iv)), and since $\alpha_0$ is
topologically nilpotent, for every open ideal $I$ of
$W(\bE)$ there exists $n\in\N$ such that
$W((\alpha^n_0\cK)^c)\subset I$. Consequently, the
series $\sum_{n\in\N}y_n$ converges to an element
$y\in W(\cK)$, and $\lim_{n\to+\infty}x_n=0$, so finally
$x=\underline\alpha\cdot y$, which concludes the proof.
\end{proof}

\begin{corollary}\label{cor_perfectoid-quot}
Let $A$ be any perfectoid ring, $\cK\subset\bE:=\bE(A)$
an ideal such that $\Phi_\bE(\cK)=\cK$, and endow $A/\{\cK\}$
(resp. $\bE/\!\cK^c$) with the quotient topology induced
via the projection $\pi_A:A\to A/\{\cK\}$ (resp.
$\pi_\bE:\bE\to\bE/\!\cK^c$). Then $A/\{\cK\}$ is
perfectoid, and there exists a natural isomorphism of
topological rings
$$
\omega:\bE(A/\{\cK\})\isom\bE/\!\cK^c
\qquad\text{such that}\qquad
\omega\circ\bE(\pi_A)=\pi_\bE.
$$
\end{corollary}
\begin{proof} The topological $\F_p$-algebra $\bE/\cK^c$
is complete and separated, and clearly $\Phi_\bE(\cK^c)=\cK^c$,
so $\bE/\cK^c$ is perfectoid, and we have
$$
W(\cK^c)=\Ker\,(W(\pi_\bE):W(\bE)\to W(\bE/\!\cK^c)).
$$
It follows that the kernel of the induced projection
$\pi'_A:A\to A\otimes_{W(\bE)}W(\bE/\!\cK^c)$ is naturally
identified with $u_A(W(\cK^c))$. Moreover, since
$W(\bE/\!\cK^c)$ is separated (lemma \ref{lem_Witt-limit}(ii)),
$W(\cK^c)$ is a closed subset of $W(\bE)$, and taking
into account theorem \ref{th_taut-two}(i) we deduce that
$\Ker\,\pi'_A=\{\cK\}$, so we get a cocartesian diagram
$$
\cD
\qquad :\qquad
{\diagram W(\bE) \ar[rr]^-{W(\pi_\bE)}
\ar[d]_{u_A} & & W(\bE/\!\cK) \ar[d]^{u'} \\
A \ar[rr]^-{\pi_A} & & A/\{\cK\}
\enddiagram}$$
and notice that both $u_A$ and $\pi_A$ are open and
surjective (lemma \ref{lem_was-third-cond}(i)), so the
same holds for $u'$, and therefore $A/\{\cK\}$ is perfectoid
(proposition \ref{prop_bar-W-perfect}(i)). Let also
$\pi:A\to A/pA$ be the proejction; taking into
account remark \ref{rem_general}(ii), there follows a
commutative diagram
$$
\bE(\cD\otimes_\Z\F_p)
\qquad :\qquad
{\diagram \bE \ar[rrr]^-{\pi_\bE} \ar[d]_{\bE(\pi)} & & &
\bE/\!\cK \ar[d]^{\bE(u')} \\
\bE(A/pA) \ar[rrr]^-{\bE(\pi_A\otimes_\Z\F_p)} & & &
\bE(A/(\{\cK\}+pA))
\enddiagram}$$
whose vertical arrows are both isomorphisms. The assertion
follows easily.
\end{proof}

\begin{remark}\label{rem_perfectoid-quot}
(i)\ \
In the situation of corollary \ref{cor_perfectoid-quot},
let $(\alpha_n~|~n\in\N)$ be any distinguished element of
$\Ker\,u_A$, and suppose that $\cK\subset\Ann_\bE(\alpha_0)$.
Then the inclusion $0\subset\cK^c$ is taut, and therefore
theorem \ref{th_taut-one}(ii) yields a natural
$W(\bE)$-linear identification
$$
\cK^c\isom\{\cK\}
\qquad
\beta\mapsto\bar u_A(\beta).
$$
Especially, $\{\cK\}\subset\Ann_A(p)$.

(ii)\ \
Let $\cI,\cK\subset\bE$ be two ideals such the topology
of $\bE$ agrees with the $\cI$-adic topology, and suppose
that
$$
\cK\cI=0
\qquad\text{and}\qquad
\Phi_\bE(\cK)=\cK.
$$
Notice that $\cK=\cK^{\La 1\Ra}$, and $\alpha_0^n\in\cI$
for some $n\in\N$, from which it follows easily that
$\cK\subset\Ann_\bE(\alpha_0)$. Moreover, $\cK\cap\cI=0$;
indeed, if $x\in\cK\cap\cI$, then $x^2=0$, hence $x=0$.
Thus, $\cK+\cI^n=\cK\oplus\cI^n$ for every $n\in\N$, and
therefore
$$
\bigcap_{n\in\N}(\cK+\cI^n)=\cK\oplus\bigcap_{n\in\N}\cI^n=\cK
$$
{\em i.e.} $\cK$ is closed in the topology of $\bE$,
so corollary \ref{cor_perfectoid-quot} says that $\bE/\cK$
and $A/\{\cK\}$ are perfectoid for their quotient topologies,
and there exists an isomorphism of topological rings
$$
\omega:\bE(A/\{\cK\})\isom\bE/\!\cK
\qquad\text{such that}\qquad
\omega\circ\bE(\pi_A)=\pi_\bE
$$
where $\pi_A:A\to A/\{\cK\}$ and $\pi_\bE:\bE\to\bE/\!\cK$
are the projections.
\end{remark}

\begin{corollary}\label{cor_taut-two}
Let $A$ be any perfectoid ring, and $a_\bullet:=(a_1,\dots,a_k)$
a finite system of elements of\/ $\bE:=\bE(A)$. Denote by
$J\subset\bE$ (resp. $\cJ\subset W(\bE)$) the ideal generated
by $a_\bullet$ (resp. by $\tau_\bE(a_1),\dots,\tau_\bE(a_k)$),
and define the ideal $[a_\bullet]^{\La\lambda\Ra}\subset W(\bE)$
as in remark {\em\ref{rem_Witt-are-f-adic}(iii)}, for every
$\lambda\in\N[1/p]$. We regard $A$ as a $W(\bE)$-algebra,
via the map $u_A$. Then we have :
\begin{enumerate}
\item
$[a_\bullet]^{\La\lambda\Ra}A\subset\{J^{\La\lambda\Ra}\bE\}
\subset[a_\bullet]^{\La\lambda'\Ra}A$ \ \
for every $\lambda'<\lambda$ in $\N[1/p]$.
\item
The following conditions are equivalent :
\begin{enumerate}
\item
$J$ is open in $\bE$
\item
$\cJ A$ is open in $A$
\item
$[a_\bullet]^{\La\lambda\Ra}A$ is open in $A$ for
every $\lambda\in\N[1/p]$
\item
There exists a strictly positive $\lambda\in\N[1/p]$
such that $[a_\bullet]^{\La\lambda\Ra}A$ is open in $A$
\end{enumerate}
and if these conditions hold, then
$[a_\bullet]^{\La\lambda\Ra}A=\{J^{\La\lambda\Ra}\bE\}$
\ \ for every $\lambda\in\N[1/p]$.
\item
If $I$ is any ideal of definition of\/ $\bE$, the following
holds :
\begin{enumerate}
\item
$\{I^{\La\lambda\Ra}\bE\}$ is a topologically nilpotent
open ideal of\/ $A$ for every $\lambda>0$ in $\N[1/p]$.
\item
$\bigcap_{n\in\N}\{J^{\La\lambda\Ra}\bE+I^{\La n\Ra}\bE\}\subset
\{J^{\La\lambda'\Ra}\bE\}$
\ \ for every $\lambda'<\lambda$ in $\N[1/p]$.
\end{enumerate}
\end{enumerate}
\end{corollary}
\begin{proof}(i): The first inclusion is clear. For
the second, notice that
$\{J^{\La\lambda\Ra}\bE\}=W(J^{\La\lambda\Ra})^c\cdot A$
by theorem \ref{th_taut-two}(i) and remark
\ref{rem_beta-taut}(iii,v). On the other hand, pick any
$\lambda''\in\N[1/p]$ such that $\lambda'<\lambda''<\lambda$;
we have
$W(J^{\La\lambda\Ra})^c=\prod_{n\in\N}(J^{\La p^n\lambda\Ra}\bE)^c$
(see remark \ref{rem_Witt-limit}(iv) and lemma
\ref{lem_mon-fract-powers}(ii.a)), and moreover
$(J^{\La p^n\lambda\Ra}\bE)^c\subset J^{\La p^n\lambda''\Ra}\bE$
for every $n\in\N$ (claim \ref{cl_B-is-A-cplt}), hence
$$
W(J^{\La\lambda\Ra})^c
\subset W(J^{\La\lambda''\Ra})
\subset[a_\bullet]^{\La\lambda'\Ra}
$$
where the last inclusion follows from proposition
\ref{prop_morel}(i).

(ii): Taking into account lemma
\ref{lem_mon-fract-powers}(ii.a,iv), it is easily seen
that (ii.b)$\Leftrightarrow$(ii.c)$\Leftrightarrow$(ii.d).
Combining with (i), we conclude that (ii.b) holds
if and only if $\{J^{\La 1\Ra}\bE\}$ is open in $A$.
However, from theorem \ref{th_taut-two}(ii) and
remark \ref{rem_beta-taut}(iii) we see as well that
$\{J^{\La 1\Ra}\bE\}$ is open in $A$ if and only if
$J^{\La 1\Ra}\bE$ is open in $\bE$. Lastly, $J^{\La 1\Ra}\bE$
is open in $\bE$ if and only if the same holds for
$J$, due to lemma \ref{lem_mon-fract-powers}(ii.a,iv).
Thus (ii.a)$\Leftrightarrow$(ii.b). Next,
we remark :

\begin{claim}\label{claim_u-is-dense}
For every $\lambda\in\N[1/p]$ and every $n\in\N$,
the ideal generated by the system
$(\bar u_A(x)~|~x\in J^{\La\lambda\Ra}\bE)$ is
contained in $p^nA+[a_\bullet]^{\La\lambda\Ra}A$.
\end{claim}
\begin{pfclaim} The claim follows easily from proposition
\ref{prop_combinatorial} and lemma \ref{lem_drop-conditions}(iv)
(details left to the reader).
\end{pfclaim}

Now, if $\cJ A$ is open in $A$, the same holds for
$[a_\bullet]^{\La\lambda\Ra}A$, for every
$\lambda\in\N[1/p]$, in which case the latter contains
$p^n$ for some sufficiently large $n\in\N$; taking into
account claim \ref{claim_u-is-dense}, we get the second
assertion of (ii).

(iii.a): It follows easily from (ii) that $\{I^{\La\lambda\Ra}\}$
is open, and combining lemma \ref{lem_mon-fract-powers}(ii.a,iv)
with (i), it is easily seen that $\{I^{\La\lambda\Ra}\bE\}$
is topologically nilpotent : details left to the reader.

(iii.b): Fix $N\in\N$ such that $\lambda'<(1-p^{-N})\lambda$,
let $x\in J^{\La\lambda\Ra}\bE$, $y\in I^{\La n\Ra}$ be
any two elements, and $\alpha\in\N[1/p]$ any rational
number $\leq 1$; notice that
$z_\alpha:=\tau_A(x^{1-\alpha}\cdot y^\alpha)$ lies in
$\{I^{\La n/p^N\Ra}\bE\}$ if $\alpha\geq p^{-N}$, and lies
in $\{J^{\La\lambda'\Ra}\bE\}$ otherwise. However, proposition
\ref{prop_combinatorial} says that $\tau_A(x+y)$ can be
written as a $p$-adically convergent series whose terms
are of the form $u_\alpha\cdot z_\alpha$ for certain
$u_\alpha\in W(\bE)$; since (ii) implies that
$\{I^{\La n/p^N\Ra}\bE\}$ is an open ideal in $A$, we deduce
$$
\bigcap_{n\in\N}\{J^{\La\lambda\Ra}\bE+I^{\La n\Ra}\bE\}\subset
\bigcap_{n\in\N}(\{J^{\La\lambda'\Ra}\bE\}
+\{I^{\La n/p^N\Ra}\bE\})=\{J^{\La\lambda'\Ra}\bE\}
$$
where the last equality holds because $\{J^{\La\lambda'\Ra}\bE\}$
is a closed ideal, taking into account (iii.a).
\end{proof}

\begin{theorem}\label{th_adic-to-adic}
Let $f:A\to A'$ be a continuous ring homomorphism of
perfectoid rings; set $C:=\Img\,f$ and $D:=\Img\,\bE(f)$.
The following holds :

{\em (i)}\ \
$f$ is surjective (resp. bijective, resp. adic, resp. open,
resp. open and injective) if and only if the same holds for
$\bE(f)$.

{\em (ii)}\ \
$C$ is open in $A'$ if and only if $D$ is open in $\bE(A')$.

{\em (iii)}\ \
Suppose that $C$ is open in $A'$, then we have :
\begin{enumerate}
\alphaenu
\item
$A'^{\circ\circ}\subset C$.
\item
$C$ is perfectoid both for the quotient topology induced by
the projection $A\to C$, and for the subspace topology induced
by the inclusion map $C\to A'$.
\end{enumerate}
\end{theorem}
\begin{proof} Set $\bE:=\bE(A)$ and $\bE':=\bE(A')$. Also,
denote by $\cT_D$ and $\cT'_D$ (resp. $\cT_C$ and $\cT'_C$)
the topologies on $D$ (resp. on $C$) induced respectively
by the projection $\bE\to D$ and the open inclusion $D\to\bE'$
(resp. induced by the projection $A\to C$ and the open
inclusion $C\to A'$).

Let $\alpha_\bullet:=(\alpha_n~|~n\in\N)$ be a distinguished
element in the kernel of $u_A:W(\bE)\to A$, and recall that
the image of $\alpha_\bullet$ in $W(\bE')$ is a distinguished
element in the kernel of $u_{A'}$.
We denote by $\Fil^\bullet A$ and $\Fil^\bullet A'$ the
$p$-adic filtrations on $A$ and respectively $A'$;
likewise, let $\Fil^\bullet\bE$ and $\Fil^\bullet\bE'$ be the
$\alpha_0$-adic filtrations on $\bE$ and respectively $\bE'$.
Moreover, we let $A^\bullet$ be the complex concentrated
in degrees $0$ and $1$, with $A^0:=A$, $A^1:=A'$, and
$d^0:=f$. Likewise, we define the complex
$\bE^\bullet:=(\bE\xrightarrow{\bE(f)}\bE')$ in degrees
$0$ and $1$. With the foregoing filtrations, $A^\bullet$
and $\bE^\bullet$ are filtered complexes of $\Z$-modules,
and we wish to consider the associated spectral sequences,
as described in \eqref{subsec_fil-spec-seq}.
Explicitly, we have $Z(A^\bullet)^{pq}_r=0$ whenever
$p+q\neq 0,1$, and :
$$
\begin{aligned}
Z(A^\bullet)^{p,-p}_r:=&\,\{a\in\Fil^pA~|~f(a)\in\Fil^{p+r}A'\} \\
Z(A^\bullet)^{p,-p+1}_r:=&\,\Fil^pA'
\end{aligned}
\qquad
\text{for every $p\in\Z$ and $r\in\N$}
$$
and correspondingly for $Z(\bE^\bullet)^{pq}_r$. Since
$\bar u_A(\Fil^p\bE)\subset\Fil^pA$ for every $p\in\Z$,
and likewise for $\bar u_{A'}$, it is then clear that
$\bar u_A$ and $\bar u_{A'}$ restrict to well defined maps
$$
\bar u{}^{pq}_r:Z(\bE^\bullet)^{pq}_r\to Z(A^\bullet)^{pq}_r
\qquad
\text{for every $p,q\in\Z$ and $r\in\N$}.
$$
Moreover, we have $B(A^\bullet)^{pq}_r=0$ whenever
$p+q\neq 0,1$, and :
$$
\begin{aligned}
B(A^\bullet)^{p,-p}_r:=&\,\Fil^{p+1}A\cap Z(A^\bullet)^{p,-p}_r \\
B(A^\bullet)^{p,-p+1}_r:=&\,
(\Fil^{p+1}A'+f(\Fil^{p-r+1}A))\cap\Fil^pA'
\end{aligned}
\qquad
\text{for every $p\in\Z$ and $r\in\N$}
$$
and correspondingly for $B(\bE^\bullet)^{pq}_r$.

\begin{claim}\label{cl_dylan-come-back}
$\bar u{}^{a,b}_r(B(\bE^\bullet)^{a,b}_r)\subset B(A^\bullet)^{a,b}_r$
for every $a,b\in\Z$ and $r\in\N$.
\end{claim}
\begin{pfclaim} The assertion is trivial for $a+b\neq 0,1$.
If $b=-a$, the assertion follows directly from our
explicit descriptions of $B(\bE^\bullet)^{a,-a}_r$ and
$B(A^\bullet)^{a,-a}_r$. Next, notice that
$B(\bE^\bullet)^{a,-a+1}_0=\Fil^{a+1}\bE'$, and likewise for
$B(A^\bullet)^{a,-a+1}_0$; the assertion follows then when
$r=0$ and $b=-a+1$. For the case where $r>0$ and  $b=-a+1$,
lemma \ref{lem_notice-that} shows that
$$
B(\bE^\bullet)^{a,-a+1}_r=
\Fil^{a+1}\bE'+f(Z(\bE^\bullet)^{a-r+1,-a+r-1}_{r-1})
$$
and likewise for $B(A^\bullet)^{a,-a+1}_r$. Hence, let
$x\in\Fil^{p+1}\bE'$ and $y\in Z(\bE^\bullet)^{a-r+1,-a+r-1}_{r-1}$
and set $z:=x+\bE(f)(y)$; we need to check that
$\bar u_{A'}(z)\in\Fil^{a+1}A'+f(Z(A^\bullet)^{a-r+1,-a+r-1}_{r-1})$.
However, proposition \ref{prop_combinatorial} says that
$\bar u_{A'}(z)=\bar u_{A'}(x)+\bar u_{A'}\circ\bE(f)(y)+w$,
where $w$ is the limit of a $p$-adically convergent series
whose terms are $\Z_p$-linear combinations of products of
the form $p^n\cdot\bar u_{A'}(x^{\sigma_0}\cdot\bE(f)(y)^{\sigma_1})$,
with $n\in\N\setminus\{0\}$, $\sigma_0,\sigma_1\in\Q_+$, and
$\sigma_0+\sigma_1=1$. Now :
$$
\bar u_{A'}\circ\bE(f)(y)=f(\bar u_A(y))\in
f(Z(A^\bullet)^{a-r+1,-a+r-1}_{r-1})
\qquad\text{and}\qquad
\bar u_{A'}(x)\in\Fil^{a+1}A'.
$$
Moreover, since $f(\bar u_A(y))\in\Fil^aA'$, we see that
$p^m$ divides $\bar u_{A'}(x^{\sigma_0}\cdot\bE(f)(y)^{\sigma_1})=
\bar u_{A'}(x)^{\sigma_0}\cdot f(\bar u_A(y))^{\sigma_1}$ in $A'$,
for every integer
$m\leq\sigma_0\cdot(a+1)+\sigma_1\cdot a=a+\sigma_0$. Thus,
$w\in\Fil^{a+1}A'$, whence the contention.
\end{pfclaim}

Next, notice that $E(A^\bullet)^{p,-p}_0=\gr^pA$ and
$E(A^\bullet)^{p,-p+1}_0=\gr^pA'$ for every $p\in\Z$,
and likewise for $E(\bE^\bullet)^{pq}_0$, whenever $p+q=0,1$.
Especially, $E(A^\bullet)^{a,b}_r$ is a $\Z/p\Z$-module
for every $a,b\in\Z$ and $r\in\N$. Invoking proposition
\ref{prop_combinatorial} again, we deduce that the composition
$$
Z(\bE^\bullet)^{pq}_r\xrightarrow{\bar u{}^{pq}_r}
Z(A^\bullet)^{pq}_r\to E(A^\bullet)^{pq}_r
$$
is an additive map, and combining with claim
\ref{cl_dylan-come-back}, it follows that the latter
factors through a well defined group homomorphism
$$
v^{pq}_r:E(\bE^\bullet)^{pq}_r\to E(A^\bullet)^{pq}_r
\qquad
\text{for every $p,q\in\Z$ and $r\in\N$}.
$$
Lastly, a simple inspection shows that the system of
maps $v^{\bullet\bullet}_\bullet$ yields a morphism of spectral
sequences $E(\bE^\bullet)^{\bullet\bullet}_\bullet\to
E(A^\bullet)^{\bullet\bullet}_\bullet$. Since $v^{pq}_0$ is
an isomorphism for every $p,q\in\Z$, it follows that
$v^{\bullet\bullet}_\bullet$ is an isomorphism of spectral sequences.

Let us endow $\Coker\,f$ with the filtration induced by
$\Fil^\bullet A'$ : namely
$$
\Fil^p(\Coker\,f):=\Img(\Fil^pA'\to\Coker\,f)
\qquad
\text{for every $p\in\Z$}.
$$
Likewise, we define a filtration $\Fil^\bullet\Coker\,\bE(f)$
on $\Coker\,\bE(f)$. Let moreover $\cT_{A',p}$ be the $p$-adic
topology on $A'$, and $\cT_{\bE',\alpha_0}$ the $\alpha_0$-adic
topology on $\bE'$.

\begin{claim}\label{cl_condition-for-open}
$C$ is open in $(A',\cT_{A',p})$ (resp. $D$ is open in
$(\bE',\cT_{\bE',\alpha_0})$) if and only if there
exists $i\in\N$ such that $\gr^i\Coker\,f=0$ (resp. such
that $\gr^i\Coker\,\bE(f)=0$).
\end{claim}
\begin{pfclaim} Clearly, if $C$ is open in $(A',\cT_{A',p})$,
then we find $i\in\N$ with $\gr^i\Coker\,f=0$. Conversely, if
the latter condition holds, set $M:=\Fil^iA'+C\subset A'$; we
deduce that $M=pM+C$, whence $M=C$, by \cite[Th.8.4]{Mat} and
lemma \ref{lem_perfectoid}(v); {\em i.e.} $\Fil^iA'\subset C$,
so $C$ is open in $(A',\cT_{A',p})$. Since $\bE$ and $\bE'$
are complete and separated for their $\alpha_0$-adic
topologies (see \eqref{subsec_equivalence}), the same
argument shows the assertion for $D$.
\end{pfclaim}

According to \eqref{subsec_Fil-toSpSeq}, the spectral
sequences $E(A^\bullet)^{\bullet\bullet}_\bullet$ and
$E(\bE^\bullet)^{\bullet\bullet}_\bullet$ admit natural
abutments, and by inspecting the constructions, we deduce
natural isomorphisms :
$$
E(A^\bullet)^{a,1-a}_\infty\isom\gr^a\Coker\,f
\qquad
E(\bE^\bullet)^{a,1-a}_\infty\isom\gr^a\Coker\,\bE(f)
\qquad
\text{for every $a\in\Z$}.
$$
Moreover, by proposition \ref{prop_convergence-simple}
they converge in degree $1$. Since these spectral sequences
are isomorphic, combining with claim
\ref{cl_condition-for-open} we conclude that $C$ is open in
$(A',\cT_{A',p})$ if and only if $D$ is open in
$(\bE',\cT_{\bE',\alpha_0})$.

(ii): By the foregoing, in order to prove the assertion
we may assume that $C$ is open in $(A',\cT_{A',p})$ and
$D$ is open in $(\bE',\cT_{\alpha_0,\bE'})$. Hence,
$\alpha_0^n\bE'\subset D$ for some $n\in\N$, and since
$D$ and $\bE'$ are perfect rings, it follows easily that
$\alpha_0\bE'\subset D$, {\em i.e.} $\gr^1\Coker\,\bE(f)=0$.
By the foregoing, we deduce that $\gr^1\Coker\,f=0$ as
well, {\em i.e.} $pA'\subset C$. Endow $\bE'/\alpha_0\bE'$
and $A'/pA'$ with the quotient topologies $\cT_{\bE'/\alpha_0\bE'}$
and $\cT_{A'/pA'}$ induced by $\bE'$ and $A'$, and recall that
$\bar u_{A'/pA'}$ induces an isomorphism of topological rings
$$
\omega:(\bE'/\alpha_0\bE',\cT_{\bE'/\alpha_0\bE'})\isom
(A'/pA',\cT_{A'/pA'})
$$
(remark \ref{rem_nice-topology}(ii)); clearly
$\omega(C/pA')=D/\alpha_0\bE'$, and especially $C/pA'$
is open in $A'/pA'$ if and only if $D/\alpha_0\bE'$ is
open in $\bE'/\alpha_0\bE'$. The contention follows
immediately.

(iii.a): In view of (ii), the subring $D$ is open in $\bE'$;
since $D$ is perfect, it follows easily that
$\bE'^{\circ\circ}\subset D$. Endow $D/\alpha_0\bE'$ (resp.
$C/pA'$) with the topology $\cT'_{D/\alpha_0\bE'}$ (resp.
$\cT'_{C/pA'}$) induced by $\cT'_D$ via the projection
$D\to D/\alpha_0\bE'$ (resp. by $\cT'_C$ via the projection
$C\to C/pA'$). Then we have
$(D,\cT'_D)^{\circ\circ}=\bE'^{\circ\circ}$, and notice that
$$
(D/\alpha_0\bE',\cT'_{D/\alpha_0\bE'})^{\circ\circ}=
(D,\cT'_D)^{\circ\circ}/\alpha_0\bE'
\qquad\text{and}\qquad
(C/pA',\cT'_{C/pA'})^{\circ\circ}=(C,\cT'_C)^{\circ\circ}/pA'.
$$
Moreover, according to lemma \ref{lem_no-need-of-dense}(i),
the topology $\cT'_{D/\alpha_0\bE'}$ agrees with the topology
induced by $\cT_{\bE'/\alpha_0\bE'}$ via the inclusion map
$D/\alpha_0\bE'\to\bE'/\alpha_0\bE'$, and likewise for
$\cT'_{C/pA'}$; it follows that the foregoing isomorphism
$\omega$ restricts to an isomorphism of topological rings
$$
(C/pA',\cT'_{C/pA'})\isom(D/\alpha_0\bE',\cT'_{D/\alpha_0\bE'}).
$$
Since $(D/\alpha_0\bE',\cT'_{D/\alpha_0\bE'})^{\circ\circ}=
\bE'^{\circ\circ}/\alpha_0\bE'=
(\bE'/\alpha_0,\cT_{\bE'/\alpha_0})^{\circ\circ}$, we get
$(C/pA',\cT'_{C/pA'})^{\circ\circ}=(A'/pA',\cT_{A'/pA'})^{\circ\circ}=
A'^{\circ\circ}/pA'$, and finally $A'^{\circ\circ}=(C,\cT'_C)^{\circ\circ}$.

(i): If $f$ is bijective, then clearly the same holds for
$\bE(f)$; conversely, if $\bE(f)$ is bijective, then
the same holds for
$$
W(\bE)/\alpha_\bullet W(\bE)\otimes_{W(\bE)}W(\bE(f)):
W(\bE)/\alpha_\bullet W(\bE)\to W(\bE')/\alpha_\bullet W(\bE').
$$
But the image of $\alpha_\bullet$ is still distinguished
in $W(\bE')$, and lies in $\Ker\,u_{A'}$, hence the latter
map is naturally identified with $f$; especially, $f$
is bijective.

$\bullet$\ \
If either $f$ or $\bE(f)$ is surjective, then
$pA'\subset C$ and $\alpha_0\bE'\subset D$ by the
proof of (ii), and moreover the isomorphism
$\bE'/\alpha_0\bE'\isom A'/pA'$ induced by $\bar u_{A/pA}$
restricts to an isomorphism $D/\alpha_0\bE'\isom C/pA'$. Then
clearly $C=A'$ if and only if $D=\bE'$.

$\bullet$\ \
Next, let $\cJ\subset\bE$ be any ideal of definition,
$(\beta_1,\dots,\beta_k)$ a finite system of generators
for $\cJ$, and $J\subset A$ the ideal generated by
$(\bar u_A(\beta_1),\dots,\bar u_A(\beta_k))$. By
corollary \ref{cor_taut-two}(ii), the ideal $JA'$ is open
in $A'$ if and only if $\cJ\bE'$ is open in $\bE'$. Moreover,
$\cJ\bE'$ (resp. $JA'$) is topologically nilpotent in $\bE'$
(resp. in $A'$), since the same holds for $\cJ$ (resp.
for $J$); thus, $\cJ\bE'$ is an ideal of adic definition
for $\bE'$ if and only if $JA'$ is an ideal of adic
definition for $A'$, {\em i.e.} $f$ is adic if and only if
$\bE(f)$ is adic.

$\bullet$\ \
If either $f$ of $\bE(f)$ is open, we know already that $C$
is open in $A'$ and $D$ is open in $\bE'$, by (ii). Also, by
the foregoing we know already that $f$ is adic if and only if
the same holds for $\bE(f)$; combining with proposition
\ref{prop_f-adics}(v), we conclude that $f$ is open if and
only if the same holds for $\bE(f)$.

$\bullet$\ \
Lastly, if either $f$ or $\bE(f)$ is open and injective,
by the foregoing we know already that both $f$ and $\bE(f)$
are open, and then $A'^{\circ\circ}\subset C$ and
$\bE'^{\circ\circ}\subset D$, by (iii.a).
Endow $\bar A{}':=A'/A'^{\circ\circ}$,
$\bar\bE{}':=\bE'/\bE'^{\circ\circ}$, $\bar C:=C/A'^{\circ\circ}$
and $\bar D:=D/\bE'^{\circ\circ}$ with their discrete topologies;
then the foregoing isomorphism $\omega$ induces a ring
isomorphism $\bar\bE{}'\isom\bar A{}'$, restricting to
a ring isomorphism $\bar D\isom\bar C$. Moreover, we get
natural isomorphisms of topological rings :
$$
(C,\cT'_C)\isom\bar C\times_{\bar A{}'}A'
\qquad\text{and}\qquad
(D,\cT'_D)\isom\bar D\times_{\bar\bE{}'}\bE'.
$$
Taking into account proposition \ref{prop_construct-new-perfs},
it follows that the projection $A\to(C,\cT'_C)$ is an isomorphism
of topological rings if and only if the same holds for the
projection $\bE\to(D,\cT'_D)$; {\em i.e.} $f$ is open and
injective if and only if the same holds for $\bE(f)$.

(iii.b): The map $\bE(f)$ is a composition of continuous ring
homomorphisms :
$$
\bE\xrightarrow{g_1}(D,\cT_D)\xrightarrow{g_2}(D,\cT'_D)
\xrightarrow{g_3}\bE'
$$
where $g_1$ is open and surjective, $g_2$ is bijective, and
$g_3$ is an open injective map. Since $\Ker\,g_1=\Ker\,\bE(f)$
is a closed ideal of $\bE$, we see that $(D,\cT_D)$ is a
perfect, separated and complete topological ring (proposition
\ref{prop_replaces-Mat-Th.8.1}(v)) whose topology is $I$-adic,
for an ideal $I\subset D$ of finite type, {\em i.e.} $(D,\cT_D)$
is perfectoid, and the same holds for $(D,\cT'_D)$, by virtue
of corollaries \ref{cor_Tate}(i) and \ref{cor_f-adics}(iii).
Since the category of perfectoid $A$-algebras is equivalent to
that perfectoid $\bE$-algebras (remark \ref{rem_nice-topology}(i))
we deduce a corresponding factorization of $f$ :
$$
A\xrightarrow{f_1}C_1\xrightarrow{f_2}C_2\xrightarrow{f_3}A'
$$
where $C_1$ and $C_2$ are perfectoid $A$-algebras with
isomorphisms  of perfectoid $\bE$-algebras
$$
\bE(C_1)\isom(D,\cT_D)
\qquad\text{and}\qquad
\bE(C_2)\isom(D,\cT'_D).
$$
Moreover, $f_1$ is open and surjective, $f_2$ is bijective,
and $f_3$ is open and injective, by (i). But then $C_1$
is isomorphic to $(C,\cT_C)$ and $C_2$ is isomorphic to
$(C,\cT'_C)$, whence the contention.
\end{proof}

In the situation of theorem \ref{th_taut-two}, let
$\cK$ and $\cK'$ be any two closed and taut ideals
of $\bE$; we would like to show that $\cK=\cK'$ if
and only if $\{\cK\}=\{\cK'\}$, and taking into account
theorem \ref{th_taut-two}(iii), it would suffice to
prove that $\{\cK\cap\cK'\}=\{\cK\}\cap\{\cK'\}$.
We do not know whether the latter identity always holds,
but we have at least the following :

\begin{theorem}\label{th_taut-three}
Let $A$ be any perfectoid ring, $\cK$ and $\cK'$ two
taut and closed ideals of\/ $\bE:=\bE(A)$, and
$(\alpha_n~|~n\in\N)$ any distinguished element of
$\Ker\,u_A$. Suppose that either one of the following
conditions holds :
\begin{enumerate}
\alphaenu
\item
$\alpha_0\in\cK'$
\item
$\cK'$ is strictly taut and $\alpha_0^n\subset\cK'$
for some sufficiently large $n\in\N$.
\end{enumerate}
Then $\{\cK\cap\cK'\}=\{\cK\}\cap\{\cK'\}$.
\end{theorem}
\begin{proof} Suppose that (a) holds; in this case,
both inclusions $\cK'\subset\bE$ and
$\cK'\cap\cK\subset\cK$ are taut, by remark
\ref{rem_beta-taut}(iv). There follows a commutative
diagram of $A/pA$-modules
$$
\xymatrix{ \cK/(\cK'\cap\cK) \ar[r] \ar[d] &
\{\cK\}/\{\cK'\cap\cK\} \ar[d] \\
\bE/\cK' \ar[r] & A/\{\cK'\}
}$$
both of whose horizontal arrows are isomorphisms, by
theorem \ref{th_taut-one}(ii), and whose left vertical
arrow is clearly injective. Then the right vertical
arrow is injective as well, whence the contention, in
this case.

In case (b) holds, pick $\eps\in\N[1/p]$ with $0<\eps<1$
and such that $\cK'$ is $\alpha_0^{1-\eps}$-taut. We set
$$
\cI_n:=\{x\in\bE~|~\alpha_0^{n\eps}x\in\cI\}
\qquad
\text{for every $n\in\N$ and every ideal $\cI\subset\bE$}.
$$

\begin{claim}\label{cl_filter-K-prime}
With the foregoing notation we have :
\begin{enumerate}
\item
$\cK'_n$ is a closed ideal of $\bE$, for every $n\in\N$.
\item
The inclusion $\cK'_n\subset\cK'_{n+1}$ is taut for
every $n\in\N$.
\item
$\cK'_0=\cK'$ and $\cK'_r=\bE$ for every sufficiently
large $r\in\N$.
\end{enumerate}
\end{claim}
\begin{pfclaim} Since $\cK'$ is closed, it is easily
seen that (i) holds. Say that $\alpha_0^n\in\cK'$; then
clearly $\cK'_k=\bE$ for every integer $k\geq \eps^{-1}n$.
Next, directly from the definition we see that, for
any two ideals $\cI,\cJ\subset\bE$ and any $n\in\N$
we have
$$
\cJ\cdot\cI_n\subset(\cJ\cI)_n
\qquad
\alpha_0^\eps\cdot\cI_{n+1}\subset\cI_n
\qquad
(\cI_n)^p\subset(\cI^p)_{np}
\qquad
(\Phi^{-1}_\bE\cI)_n=\Phi^{-1}_\bE(\cI_{np}).
$$
Therefore :
$$
\begin{aligned}
\alpha_0\cdot\Phi^{-1}_\bE(\cK'_{n+1})^p \subset\, &
\alpha_0\cdot\Phi^{-1}_\bE((\cK^{\prime p})_{(n+1)p}) \\
=\, & \alpha_0\cdot(\Phi^{-1}_\bE\cK^{\prime p})_{n+1} \\
\subset\, & \alpha_0^{1-\eps}\cdot(\Phi^{-1}_\bE\cK^{\prime p})_n \\
\subset\, & (\alpha_0^{1-\eps}\cdot\Phi^{-1}_\bE\cK^{\prime p})_n \\
\subset\, & \cK'_n
\end{aligned}
$$
as required.
\end{pfclaim}

By claim \ref{cl_filter-K-prime} we have $\cK'_r=\bE$
for some $r\in\N$, and for every $n\in\N$ we get
a commutative diagram of $A/pA$-modules
$$
\xymatrix{ \{\cK\cap\cK'_{n+1}\}/\{\cK\cap\cK'_n\}
\ar[r] \ar[d] & \{\cK'_{n+1}\}/\{\cK'_n\} \ar[d] \\
(\cK\cap\cK'_{n+1})/(\cK\cap\cK'_n) \ar[r] &
\cK'_{n+1}/\cK'_n
}$$
whose vertical arrows are isomorphisms, by theorem
\ref{th_taut-one}(ii), and whose bottom horizontal
arrow is injective. Then the top horizontal arrow is
injective as well, so that
$$
\{\cK\cap\cK'_{n+1}\}\cap\{\cK'_n\}\subset\{\cK\cap\cK'_n\}
\qquad
\text{for every $n\in\N$}.
$$
We may now prove that $\{\cK\}\cap\{\cK'_n\}=\{\cK\cap\cK'_n\}$,
by descending induction on $n$. Indeed, the assertion is
obvious in case $n\geq r$. Suppose that the identity has
already been shown for some $n>0$; then we get
$$
\{\cK\}\cap\{\cK'_{n-1}\}=\{\cK\}\cap\{\cK'_n\}\cap\{\cK'_{n-1}\}=
\{\cK\cap\cK'_n\}\cap\{\cK'_{n-1}\}\subset\{\cK\cap\cK'_{n-1}\}
$$
and the converse inclusion is obvious, so the sought
identity holds for $n-1$. Letting $n=0$, we get the
assertion.
\end{proof}

In the same vein, let us also point out :

\begin{lemma}\label{lem_product-of-taut}
In the situation of theorem {\em\ref{th_taut-three}},
let $\lambda_1,\lambda_2\in\N[1/p]$ be two rationals
such that $\lambda_1+\lambda_2\leq 1$, and
$\cJ_1,\cJ_2\subset\bE$ two ideals such that $\cJ_i$
is $\alpha_0^{\lambda_i}$-taut for $i=1,2$. Then
$$
(\{\cJ_1\}\{\cJ_2\})^c=\{\cJ_1\cJ_2\}.
$$
\end{lemma}
\begin{proof} The inclusion
$\{\cJ_1\}\{\cJ_2\}\subset\{\cJ_1\cJ_2\}$ is obvious
from the definitions, so it suffices to check that
for every $x\in\cJ_1\cJ_2$ we have
$\bar u_A(x)\in(\{\cJ_1\}\{\cJ_2\})^c$. However,
every such $x$ can be written as a finite sum
$\sum_{i=1}^ny_iz_i$, with $y_1,\dots,y_n\in\cJ_1$ and
$z_1,\dots,z_n\in\cJ_2$; by proposition
\ref{prop_combinatorial} we know that $\bar u_A(x)$
is the limit of a $p$-adically convergent series of
the form $\sum_{j\in\N}w_j$, where in turn each $w_j$ is a
finite $\Z_p$-linear combination of terms of the form
$$
w^{\underline\mu}_j:=
p^j\cdot\bar u_A((y_1z_1)^{\mu_1}\cdots(y_nz_n)^{\mu_n})
\quad
\text{for certain $\mu_1,\dots,\mu_n\in p^{-j}\N$ such that
$\sum_{i=1}^n\mu_i=1$}.
$$
Thus, we are reduced to checking that any such term
lies in $\{\cJ_1\}\{\cJ_2\}$. Now, set
$y^{\underline\mu}:=y_1^{\mu_1}\cdots y_n^{\mu_n}$, and
define likewise $z^{\underline\mu}$, for every
$\underline\mu:=(\mu_1,\dots,\mu_n)\in p^{-j}\N$;
under our assumptions we have
$p=\bar u_A(\alpha_0^{\lambda_1+\lambda_2})\cdot a$ for some
$a\in A$ (lemma \ref{lem_was-third-cond}(iii)), hence
$$
w^{\underline\mu}_j=a^j\cdot
\bar u_A(\alpha_0^{j\lambda_1}\cdot y^{\underline\mu})
\cdot
\bar u_A(\alpha_0^{j\lambda_2}\cdot z^{\underline\mu})
$$
and notice that
$\alpha_0^{j\lambda_1}\cdot y^{\underline\mu}\in
\alpha_0^{j\lambda_1}\cdot\Phi_\bE^{-j}(\cJ_1^{p^j})
\subset\cJ_1$ by \eqref{eq_iterate-beta}. Likewise,
$\alpha_0^{j\lambda_2}\cdot z^{\underline\mu}\in\cJ_2$, and
the assertion follows.
\end{proof}

\begin{corollary}\label{cor_taut-three}
Let $A$ be a perfectoid ring, $\cK,\cK'\subset\bE(A)$
two taut and closed ideals, and suppose that either $\cK$
or $\cK'$ fulfills conditions {\em(a)} and {\em(b)} of
theorem {\em\ref{th_taut-three}}. Then we have
$\cK\subset\cK'$ if and only if\/ $\{\cK\}\subset\{\cK'\}$.
\end{corollary}
\begin{proof} As already announced, this follows
straightforwardly from remark \ref{rem_beta-taut}(iv) and
theorems \ref{th_taut-three} and \ref{th_taut-two}(iii).
\end{proof}

\begin{corollary}\label{cor_bingo}
In the situation of corollary {\em\ref{cor_taut-two}},
let $a'_\bullet:=(a'_1,\dots,a'_t)$ be another finite
system of elements of\/ $\bE$, and $J'\subset\bE$ the
ideal generated by $a'_\bullet$. Define the ideal
$[a'_\bullet]^{\La\beta\Ra}\subset W(\bE)$ as in remark
{\em\ref{rem_Witt-are-f-adic}(iii)}, for every
$\beta\in\N[1/p]$. Moreover, let
$\lambda,\mu\in\N[1/p]$ be any two rational numbers.
We have :
\begin{enumerate}
\item
The following conditions are equivalent :
\begin{enumerate}
\item
$J^{\La\lambda\Ra}\bE\subset J'{}^{\lfloor\mu\rfloor}\bE$.
\item
$[a_\bullet]^{\La\lambda\Ra}A\subset[a'_\bullet]^{\La\mu-\eps\Ra}A$
for every $\eps\in\N[1/p]$ with $0<\eps\leq\mu$.
\end{enumerate}
\item
Suppose that either $J$ or $J'$ is open. Then the
following conditions are equivalent :
\begin{enumerate}
\item
$J^{\La\lambda\Ra}\bE\subset J'{}^{\La\mu\Ra}\bE$.
\item
$[a_\bullet]^{\La\lambda\Ra}A\subset[a'_\bullet]^{\La\mu\Ra}A$.
\end{enumerate}
\item
For every $l,m\in\N[1/p]$ such that
$l\mu+m\lambda<\lambda\mu$ we have
$$
[a_\bullet]^{\La\lambda\Ra}A\cap[a'_\bullet]^{\La\mu\Ra}A
\subset [a_\bullet]^{\La l\Ra}\cdot[a'_\bullet]^{\La m\Ra}A.
$$
\end{enumerate}
\end{corollary}
\begin{proof}(i): Suppose that (i.a) holds; we deduce :
$$
[a_\bullet]^{\La\lambda\Ra}A\subset
\{J^{\La\lambda\Ra}\bE\}\subset\{J'{}^{\La\mu-\eps/p\Ra}\bE\}
\subset[a'_\bullet]^{\La\mu-\eps\Ra}A
\qquad
\text{for every $\eps$ as in (i.b)}
$$
where the first and last inclusions follow from
corollary \ref{cor_taut-two}(i).

Conversely, suppose that (i.b) holds; fix a finite
system $b_\bullet:=(b_1,\dots,b_r)$ of generators for
an ideal of definition $I$ of the perfectoid ring
$\bE$, and for every $n\in\N$, let $b_\bullet^n$ be
the finite system consisting of the products of the
form $b_{j(1)}\cdots b_{j(n)}$, where $j$ ranges over
all the mappings $\{1,\dots,n\}\to\{1,\dots,r\}$.
For every $n\in\N$, we may then form the system
$(a'_\bullet,b^n_\bullet)$ which is the union of the
systems $a'_\bullet$ and $b^n_\bullet$, and for every
$\lambda\in\N[1/p]$ consider the ideal
$[a'_\bullet,b^n_\bullet]^{\La\lambda\Ra}\subset W(\bE)$
attached as in remark \ref{rem_Witt-are-f-adic}(iii) to
the system $(a'_\bullet,b^n_\bullet)$. With this notation,
corollary \ref{cor_taut-two}(ii) yields :
$$
[a_\bullet]^{\La\lambda\Ra}A\subset
[a'_\bullet,b^n_\bullet]^{\La\mu-\eps/p\Ra}A=
\{(J'+I^n)^{\La\mu-\eps/p\Ra}\bE\}
\qquad
\text{for every $n\in\N$ and $\eps$ as in (i.a)}.
$$
However, since the topology of $A$ is coarser than
the $p$-adic topology, claim \ref{claim_u-is-dense}
implies that $[a_\bullet]^{\La\lambda\Ra}A$ is a
dense subset of $\{J^{\La\lambda\Ra}\bE\}$; but
$\{(J'+I^n)^{\La\mu-\eps/p\Ra}\bE\}$ is a closed ideal, so
\set\begin{equation}\label{eq_J-I-J-prime}
\{J^{\La\lambda\Ra}\bE\}\subset\{(J'+I^n)^{\La\mu-\eps/p\Ra}\bE\}
\qquad
\text{for every $n\in\N$ and $\eps$ as in (i.a)}.
\end{equation}
Now, the topological closure $J^{\La\lambda\Ra}\bE^c$ of
$J^{\La\lambda\Ra}\bE$ is (strictly) taut, and
$\{J^{\La\lambda\Ra}\bE\}=\{J^{\La\lambda\Ra}\bE^c\}$ (remark
\ref{rem_beta-taut}(v)); moreover, $(J'+I^n)^{\La\mu-\eps/p\Ra}\bE$
is an open ideal of $W(\bE)$, hence, from corollary
\ref{cor_taut-three} and \eqref{eq_J-I-J-prime} we
derive
$$
J^{\La\lambda\Ra}\bE^c\subset(J'+I^n)^{\La\mu-\eps/p\Ra}\bE
\qquad
\text{for every $n\in\N$ and $\eps$ as in (i.a)}
$$
and combining with claim \ref{cl_B-is-A-cplt} we conclude
that $J^{\La\lambda\Ra}\bE\subset J'{}^{\La\mu-\eps\Ra}\bE$
for every $\eps$ as in (i.a).

(ii): From lemma \ref{lem_mon-fract-powers}(ii.a,iv) it
follows easily that $J$ is open if and only if the
same holds for $J^{\La\lambda\Ra}\bE$, and likewise,
$J'$ is open if and only if the same holds for
$J'{}^{\La\mu\Ra}\bE$. Hence, suppose that $J$ is open;
if condition (ii.a) holds, then it follows that $J'$ is
open as well, and if (ii.b) holds, we use corollary
\ref{cor_taut-two}(ii) to deduce first that
$[a_\bullet]^{\La\lambda\Ra}A$ is open, so the same
holds for $[a'_\bullet]^{\La\mu\Ra}A$, and then
also for $J'$, again by corollary \ref{cor_taut-two}(ii).
In conclusion, for the proof of (ii) we may assume
that $J'$ is open. Now, suppose that (ii.a) holds;
then we get
$$
[a_\bullet]^{\La\lambda\Ra}A\subset
\{J^{\La\lambda\Ra}\bE\}\subset\{J'{}^{\La\mu\Ra}\bE\}
=[a'_\bullet]^{\La\mu\Ra}A
$$
where the first inclusion follows from
corollary \ref{cor_taut-two}(i), and the last
identity follows from corollary \ref{cor_taut-two}(ii).
Conversely, suppose that (ii.b) holds; since the
topology of $A$ is coarser that the $p$-adic topology,
claim \ref{claim_u-is-dense} implies that
$[a_\bullet]^{\La\lambda\Ra}A$ is a dense subset
of $\{J^{\La\lambda\Ra}\bE\}$, and since
$\{J'{}^{\La\mu\Ra}\bE\}$ is a closed ideal, corollary
\ref{cor_taut-two}(ii) yields
$\{J^{\La\lambda\Ra}\bE\}\subset\{J'{}^{\La\mu\Ra}\bE\}$.
Then, arguing as in the proof of (i), we easily deduce
that (ii.a) holds.

(iii): We may assume that $\lambda,\mu>0$, since otherwise
there is nothing to prove. Define $b_\bullet^n$ as in the
foregoing, for every $n\in\N$; combining corollary
\ref{cor_taut-two}(i) and theorem \ref{th_taut-three}, we get
$$
[a_\bullet,b_\bullet^n]^{\La\lambda\Ra}A\cap
[a'_\bullet,b_\bullet^n]^{\La\mu\Ra}A=
\{(J+I^n)^{\La\lambda\Ra}\bE\}\cap\{(J'+I^n)^{\La\mu\Ra}\bE\}=
\{(J+I^n)^{\La\lambda\Ra}\bE\cap(J'+I^n)^{\La\mu\Ra}\bE\}
$$
for every $n\in\N$. Now, pick $\alpha\in\N[1/p]$
such that $l/\lambda<1-\alpha$ and $m/\mu<\alpha$
and set $l':=(1-\alpha)\lambda$, $m':=\alpha\mu$.
If $x\in(J+I^n)^{\La\lambda\Ra}\bE\cap(J'+I^n)^{\La\mu\Ra}\bE$,
we may write
$$
x=x^\alpha\cdot x^{1-\alpha}\in
(J+I^n)^{\La l'\Ra}\cdot(J'+I^n)^{\La m'\Ra}\bE.
$$
Next, to ease notation, for any finite sequence
$\nu_\bullet:=(\nu_1,\dots,\nu_m)$ of elements of
$\N[1/p]$ and any sequence $(x_1,\dots,x_m)$ of
elements of $\bE$, set
$$
x_\bullet^{\nu_\bullet}:=x_1^{\nu_1}\cdots x_m^{\nu_m}
\qquad\text{and}\qquad
|\nu_\bullet|:=\nu_1+\cdots+\nu_m.
$$
We remark that the ideal
$(J+I^n)^{\La l'\Ra}\cdot(J'+I^n)^{\La m'\Ra}\bE$ is
generated by all products
$$
a_\bullet^{\nu_\bullet}\cdot a'{}_\bullet^{\nu'_\bullet}\cdot
b_\bullet^{\beta_\bullet}\cdot b_\bullet^{\beta'_\bullet}
\qquad\text{where}\qquad
|\nu_\bullet|+n^{-1}\cdot|\beta_\bullet|=l'
\qquad
|\nu'_\bullet|+n^{-1}\cdot|\beta'_\bullet|=m'.
$$
Fix $l'',m''\in\N[1/p]$ such that $l<l''<l'$ and
$m<m''<m'$, and set $\eps:=(l'-l'')+(m'-m'')$;
it is easily seen that each of the above products
lies either in $J^{\La l''\Ra}\cdot J'{}^{\La m''\Ra}\bE$
or else in $I^{\La n\eps\Ra}\bE$. Moreover, we can write
$l''=p^{-N}f$, $m'':=p^{-N}g$ for suitable $f,g,N\in\N$,
and summing up, we conclude that, for every $\delta<p^{-N}$,
the ideal
$[a_\bullet]^{\La\lambda\Ra}A\cap[a'_\bullet]^{\La\mu\Ra}A$
is contained in :
$$
\bigcap_{n\in\N}
\{J^{\La l''\Ra}\cdot J'{}^{\La m''\Ra}\bE+I^{\La n\eps\Ra}\bE\}=
\bigcap_{n\in\N}\{(J^fJ'{}^g)^{\La1/p^N\Ra}\bE+I^{\La n\eps\Ra}\bE\}
\subset\{(J^fJ'{}^g)^{\La\delta\Ra}\bE\}
$$
where the first equality follows from lemma
\ref{lem_mon-fract-powers}(ii.a,iii) the last inclusion
follows from corollary \ref{cor_taut-two}(iii.b). Let
us then choose $\delta,\delta'\in\N[1/p]$ with
$p^{-N}>\delta>\delta'$ and $\delta'$ close enough to
$p^{-N}$, so that $\delta'\cdot f>l$ and $\delta'\cdot g>m$;
using corollary \ref{cor_taut-two}(i) we get
$$
\{(J^fJ'{}^g)^{\La\delta\Ra}\bE\}\subset
[a^f_\bullet a'{}_\bullet^g]^{\La\delta'\Ra}A\subset
[a_\bullet]^{\La l\Ra}\cdot[a'_\bullet]^{\La m\Ra}A
$$
whence the contention.
\end{proof}

As an application, we can complement as follows
theorem \ref{th_double-completion}(i) :

\begin{proposition} In the situation of
theorem {\em\ref{th_double-completion}}, we have :
\begin{enumerate}
\item
\eqref{eq_old-acquaintance} is a cartesian diagram
of topological rings.
\item
All the topological rings appearing in
\eqref{eq_old-acquaintance} are perfectoid.
\end{enumerate}
\end{proposition}
\begin{proof} We know already that \eqref{eq_same-again}
is a cartesian diagram of topological rings. Moreover,
it is easily seen that the $I$-adic topology on $A$
agrees with the $(I+pA)$-adic topology, and likewise
for the $J$-adic and $(I+J)$-adic topologies, so
all the topological rings appearing in \eqref{eq_same-again}
except possibly for $(A',\cT')$ are perfectoid (proposition
\ref{prop_change-topol}(ii)). By the same token, the
$IJ$-adic topology on $A$ agrees with the $(IJ+pA)$-adic
topology, so the proposition will follow, once we have
shown the following

\begin{claim} The $IJ$-adic topology on $A$ agrees
with the linear topology $\cT$ defined by the cofiltered
system of ideals $(I^n\cap J^n~|~n\in\N)$.
\end{claim}
\begin{pfclaim}[] Set $I':=I+pA$ and $J':=J+pA$,
and notice that
$I^{k+N}\subset I'{}^{N+k}\subset I^{k+1}$ for every
$k\in\N$, and likewise for $J$. Therefore, $\cT$ agrees
with the topology defined by the cofiltered system of
ideals $(I'{}^n\cap J'{}^n~|~n\in\N)$, and the $IJ$-adic
topology on $A$ agrees with the $I'J'$-adic topology.
Thus, we can replace $I,J$ by $I',J'$, and assume from
start that $p\in I\cap J$. Now, let $(\alpha_n~|~n\in\N)$
be any distinguished element in $\Ker\,u_A$. Fix finite
systems $\beta_\bullet:=(\beta_1,\dots,\beta_k)$ and
$\beta'_\bullet:=(\beta'_1,\dots,\beta'_r)$ of elements
of $\bE(A)$ such that $\bar u_{A/pA}(\beta_\bullet)$
(resp. $\bar u_{A/pA}(\beta'_\bullet)$) is a system of
generators for $I/pA$ (resp. for $J/pA$). In light of
lemma \ref{lem_was-third-cond}(iii), it follows that
the sequence
$(\bar u_A(\alpha_0),\bar u_A(\beta_1),\dots,\bar u_A(\beta_k))$
is a system of generators for $I$, and
$(\bar u_A(\alpha_0),\bar u_A(\beta'_1),\dots,
\bar u_A(\beta'_r))$ is a system of generators for
$J$. On the other hand, lemma \ref{lem_mon-fract-powers}(iv)
implies that the $I$-adic (resp. $J$-adic) topology on
$A$ agrees with the linear topology defined by the
cofiltered system of ideals
$([\alpha_0,\beta_\bullet]^{\La\lambda\Ra}A~|~\lambda\in\N[1/p])$
(resp.
$([\alpha_0,\beta'_\bullet]^{\La\lambda\Ra}A~|~\lambda\in\N[1/p])$).
Then, the claim follows from corollary \ref{cor_bingo}(iii).
\end{pfclaim}
\end{proof}

\sset\subsubsection{}\label{subsec_setup-critperf}
Let $A$ be a ring, and $a_\bullet:=(a_1,\dots,a_k)$ a finite
system of elements of $A$; denote by $I$ the ideal generated
by $a_\bullet$, and suppose that
\set\begin{equation}\label{eq_p-lies-deep}
p\in I^t
\qquad
\text{where $t:=k(p-1)+1$}.
\end{equation}
Set $J:=I^{(p)}$, the ideal generated by the
system $(a_1^p,\dots,a_k^p)$, and notice that
\set\begin{equation}\label{eq_standard-estimate}
I^t\subset J
\end{equation}
therefore $J$ is well defined independently of the
choice of $a_\bullet$ (lemma \ref{lem_special-ideals}).
Denote by $\gr^\bullet_IA$ (resp. $\gr^\bullet_JA$) the graded
ring associated with the $I$-adic (resp. $J$-adic) filtration
on $A$. Notice that both these two rings are $\F_p$-algebras.

\begin{remark}\label{rem_p-can-lie-deep}
Let $A$ be any P-ring, $I$ any ideal of definition
of $A$, and $(a_1,\dots,a_k)$ a system of generators of
$I$. By virtue of lemma \ref{lem_perfectoid}(iii) and
its proof, we may find a sequence of ideals
$(I_n~|~n\in\N)$ such that
\begin{itemize}
\item
$I_{n+1}^{(p)}=I_n$ for every $n\in\N$
\item
$I_n$ admits a system of generators consisting of
$k$ elements of $A$.
\end{itemize}
Then $I_n$ is still an ideal of definition of $A$ for
every $n\in\N$, and $p\in I^t_n$ for every sufficiently
large $n\in\N$, where $t$ is defined as in
\eqref{subsec_setup-critperf}. We may thus regard
the following result as a refinement of corollary
\ref{cor_jackob}.
\end{remark}

\begin{proposition}\label{prop_begin-criterion-perf}
In the situation of \eqref{subsec_setup-critperf}, the
following holds :
\begin{enumerate}
\item
The Frobenius endomorphism of $\gr^\bullet_IA$ factors
through a graded ring homomorphism
$$
\Phi_I:\gr^\bullet_IA\isom\gr^\bullet_JA
\qquad
(x\mod I^{n+1})\mapsto(x^p\mod J^{n+1})
\qquad
\text{for every $n\in\N$ and $x\in I^n$}.
$$
\item
Let $\cT_p$ be the $p$-adic topology on $A$, and suppose
that $(A,\cT_p)$ is perfectoid. Then $\Phi_I$ is an
isomorphism.
\end{enumerate}
\end{proposition}
\begin{proof}(i): To begin with, we remark :

\begin{claim}\label{cl_combin-difference}
For every $n\in\N$ and every $x,y\in I^n$, we have
$(x+y)^p-x^p-y^p\in J^{n+1}$.
\end{claim}
\begin{pfclaim} The difference in the claim is a sum
of terms that lie in $pI^{pn}$. Taking into account
\eqref{eq_p-lies-deep}, we see that
$pI^{pn}\subset I^{t+pn}$, so it remains only
to check that
$$
I^{t+pn}\subset J^{n+1}.
$$
We argue by induction on $n$, and notice that the
case where $n=0$ is \eqref{eq_standard-estimate}.
Hence, suppose that $r>0$, and that the sought
inclusion has already been checked for $n:=r-1$.
The ideal $I^{t+pr}$ is generated by all elements
of the type $b:=\prod_{j=1}^{t+pr}a_{\phi(j)}$, where
$\phi$ is any mapping $\{1,\dots,t+pr\}\to\{1,\dots,k\}$.
Since $r\geq 0$, it is easily seen that there exists
at least one index $i\in\{1,\dots,k\}$ such that
$\phi^{-1}(i)$ has cardinality $\geq p$. For such
index $i$, we have $b=a^p_i\cdot c$, where
$c\in I^{t+pn}$, and the inductive assumption says that
$c\in J^{n+1}$, whence $b\in J^{r+1}$, as required.
\end{pfclaim}

Let us now check that if $x\in I^n$, then $x^p\in J^n$.
Indeed, the assertion is clear if $x$ is a monomial
$\prod_{j=1}^na_{\phi(j)}$ for some mapping
$\phi:\{1,\dots,n\}\to\{1,\dots,k\}$. However, we may
write $x=x_1+\cdots+x_r$ for some $r\in\N$ and a sequence
$x_1,\dots,x_r$ of such monomials. Lastly, claim
\ref{cl_combin-difference} and an easy induction
on $r$ shows that $x^p-\sum^r_{j=1}x_j^p\in J^{n+1}$,
whence the contention. Thus, the map $\Phi_I$ is
well defined, and claim \ref{cl_combin-difference}
also shows that $\Phi_I(x+y)=\Phi_I(x)+\Phi_I(y)$ for
every $n\in\N$ and every $x,y\in\gr_I^nA$. For any
such $x$ and $y$ it is also clear that
$\Phi_I(xy)=\Phi_I(x)\cdot\Phi_I(y)$, so $\Phi_I$ is
a ring homomorphism.

(ii): Set $\bE:=\bE(A)$, and let
$\underline\alpha:=(\alpha_n~|~n\in\N)$ be a distinguished
element in $\Ker\,u_A$; arguing as in the proof of lemma
\ref{lem_was-third-cond}(i), we may assume that for every
$i=1,\dots,k$ we have $a_i=\bar u_A(\beta_i)$ for some
$\beta_i\in\bE$. Denote by $\cI\subset\bE$ the ideal generated
by $(\beta_1,\dots,\beta_k)$, and set $\cJ:=\cI^{(p)}$, the
ideal generated by $(\beta^p_1,\dots,\beta^p_k)$; from
\eqref{eq_p-lies-deep} and lemma \ref{lem_before-name}(ii)
we get
\set\begin{equation}\label{eq_same-with-alpha}
\alpha_0\in\cI^t\subset\cJ.
\end{equation}
Let $\gr^\bullet_\cI$ and $\gr^\bullet_{\!\cJ}$ be the graded
rings associated with the $\cI$-adic and $\cJ$-adic filtrations
on $\bE$.

\begin{claim}\label{cl_trickier-than-usual}
(i)\ \
For every $n\in\N$ we have :
\begin{enumerate}
\alphaenu
\item
The inclusion $\cI^{n+1}\subset\cI^n$ is $\alpha_0^{1/p}$-taut
and the inclusion $\cJ^{n+1}\subset\cJ^n$ is taut.
\item
$\{\cI^n\}=I^n$ and $\{\cJ^n\}=J^n$.
\end{enumerate}
(ii)\ \
The map $\bar u_A:\bE(A)\to A$ induces graded ring isomorphisms
$$
\gr^\bullet_\cI\isom\gr^\bullet_I
\qquad
\gr^\bullet_\cJ\isom\gr^\bullet_J.
$$
\end{claim}
\begin{pfclaim}(i.a): Clearly $\cJ=\Phi_\bE(\cI)$, therefore
it suffices to prove the assertion for $\cJ$; however
$\alpha_0\cdot\Phi^{-1}_\bE(\cJ^{np})=\alpha_0\cI^{np}
\subset\cI^{t+np}$ by \eqref{eq_same-with-alpha}, and
arguing as in the proof of claim \ref{cl_combin-difference},
we see that $\cI^{t+np}\subset\cJ^{n+1}$, whence the claim.

(i.b): We consider first the assertion for $\cI^n$. By
construction, it is clear that $a_i\in\{\cI\}$ for every
$i=1,\dots,k$, hence $I^n\subset\{\cI^n\}$.
On the other hand, since $A$ is endowed with its $p$-adic
topology, $\bE$ carries the $\alpha_0$-adic topology
(see \eqref{subsec_equivalence}), therefore $\cI^n$ is
open in $\bE$, and consequently $\{\cI^n\}$ is the ideal
generated by the system $(\bar u_A(x)~|~x\in\cI^n)$
(remark \ref{rem_beta-taut}(viii)). Thus, we come down
to checking that $\bar u_A(x)\in I^n$ for every
$x\in\cI^n$. This is clear in case $x$ is a monomial
of the form $\prod_{j=1}^n\beta_{\phi(j)}$ for some mapping
$\phi:\{1,\dots,n\}\to\{1,\dots,k\}$. In general, $x$
shall be a finite sum $x_1+\cdots+x_s$, where $x_1,\dots,x_s$
are monomials of this type. According to proposition
\ref{prop_combinatorial}, $\tau_A(x)$ can then be
expressed as a $p$-adically convergent series
$\sum_{r\in\N}p^r\cdot y_r$, where each $y_r$ is
in turn a finite sum of terms of the form $w\cdot\tau_A(z)$,
for certain $w\in\Z_p$ and $z\in\Phi^{-r}_\bE(\cI^n)^{p^r}$.
In view of (i;A) and remark \ref{rem_beta-taut}(vi), we have
$$
\alpha_0^{r/p}\cdot\Phi^{-r}_\bE(\cI^n)^{p^r}\subset\cI^{n+r}
\qquad
\text{for every $r,n\in\N$}.
$$
Pick $\pi\in A$ as in lemma \ref{lem_perfectoid}(iv); 
recalling lemma \ref{lem_was-third-cond}(iii), we deduce that
$$
p^r\cdot\bar u_A(z)\in\pi^{(p-1)r}\cdot\{\cI^{n+r}\}
\qquad
\text{for every $r\in\N$ and every $z\in\Phi^{-r}_\bE(\cI^n)^{p^r}$}.
$$
Since $\{\cI^{n+r}\}$ is closed in the $p$-adic topology
of $A$ for every $n,r\in\N$, it follows that
$\{\cI^n\}\subset I^n+\pi^{(p-1)r}\{\cI^{n+r}\}$ for
every $r\in\N$. But clearly $\pi^{(p-1)r}\in I^n$ for
every large enough $r\in\N$, so $\{\cI^n\}\subset I^n$,
as required.

To deal with $\cJ^n$, we repeat the foregoing argument,
to get : $\{\cJ^n\}\subset J^n+\{\cJ^{n+r}\}$ for every
$r\in\N$. However, clearly $\cJ^{n+r}\subset\cI^{n+r}$;
taking into account the previous case, we conclude that
$\{\cJ^n\}\subset J^n+I^s$ for every $s\in\N$. But we may
find $s\in\N$ large enough, so that $I^s\subset J^n$, so
at last $\{\cJ^n\}\subset J^n$, whence the claim.

(ii) follows immediately from (i) and theorem \ref{th_taut-one}(ii).
\end{pfclaim}

Claim \ref{cl_trickier-than-usual}(ii) yields a commutative
diagram of graded rings
$$
\xymatrix{
\gr^\bullet_\cI \ar[r] \ar[d] & \gr^\bullet_IA \ar[d]^{\Phi_I} \\
\gr^\bullet_\cJ\ar[r] & \gr^n_JA
}$$
whose horizontal arrows are isomorphisms, and whose left
vertical arrow is induced by $\Phi_\bE$. Now, since $\bE$
is perfect, the left vertical arrow is an isomorphism, so
the same holds for the right vertical arrow, as stated.
\end{proof}

\begin{corollary}\label{cor_perf-are-reduced}
Let $A$ be a perfectoid ring, and $\beta\in\bE(A)$ any
element.
\begin{enumerate}
\item
$A$ is a reduced ring.
\item
If\/ $\beta\neq 0$, then $\bar u_A(\beta)\neq 0$.
\item
If the topology of $A$ is discrete, then $A$ is an $\F_p$-algebra.
\end{enumerate}
\end{corollary}
\begin{proof}(i): Suppose that $A$ contains a nilpotent
element $x\neq 0$, and let $n>1$ be the smallest integer
such that $x^n=0$; set $y:=x^{n-1}$, so that $y\neq 0$
and $y^p=0$. By remark \ref{rem_p-can-lie-deep} we can find
an ideal of definition $I$ of $A$ fulfilling condition
\eqref{eq_p-lies-deep}; we set $J:=I^{(p)}$ and let
$\Phi_I:\gr^\bullet_IA\isom\gr^\bullet_JA$ be the isomorphism
of proposition \ref{prop_begin-criterion-perf}(ii).
Let also $r\in\N$ be the largest integer such that
$y\in I^r$, and denote by $\bar y\in\gr^r_IA$ the class
of $y$, so that $\bar y\neq 0$; but clearly
$\Phi_I(\bar y)=0$, a contradiction.

(ii) and (iii) are immediate consequences of (i).
\end{proof}

\begin{theorem}\label{th_criterium-perfect}
In the situation of \eqref{subsec_setup-critperf}, let
$\cT_I$ be the $I$-adic topology on $A$, and suppose
that $(A,\cT_I)$ is a P-ring. Define also the graded
ring homomorphism $\Phi_I$ as in proposition
{\em\ref{prop_begin-criterion-perf}(i)}. Then the
following conditions are equivalent :
\begin{enumerate}
\alphaenu
\item
$(A,\cT_I)$ is perfectoid.
\item
The map $\Phi_I$ is an isomorphism.
\end{enumerate}
\end{theorem}
\begin{proof} Suppose that (a) holds, and let $\cT_p$
be the $p$-adic topology on $A$; then $(A,\cT_p)$ is
perfectoid, by proposition \ref{prop_change-topol}(i), in
which case proposition \ref{prop_begin-criterion-perf}(ii)
says that (b) holds.

Next, suppose that (b) holds; then proposition
\ref{prop_third-equiv-cond} and example \ref{ex_perfectoid}(i)
say that $\bE:=\bE(A)$ is a perfectoid ring.
Moreover, $\bar u_{A/pA}$ is surjective (lemma
\ref{lem_was-third-cond}(i)), so for every $i=1,\dots,k$
we may find $\beta_i\in\bE$ such that
$\bar u_{A/pA}(\beta_i)$ equals the image of $a_i$ in
$A/pA$. Denote by $\cJ\subset\bE$ the ideal generated
by the system $(\beta_1,\dots,\beta_k)$; from claim
\ref{cl_drop-second} we see that $\cJ$ is an ideal of
definition of $\bE$, and $\bar u_{A/pA}$ induces a ring
isomorphism
\set\begin{equation}\label{eq_modulo-cJ}
\bE/\cJ\isom A/I.
\end{equation}
Denote by $\cJ_W\subset W(\bE)$ the ideal generated by
$(\tau_\bE(\beta_1),\dots,\tau_\bE(\beta_k))$; by
construction we have $u_A(\tau_\bE(\beta_i))-a_i\in pA$
for every $i=1,\dots,k$, and notice that $p\in I^t\subset I^2$,
so $\cJ_WA+I^2=I$ and Nakayama's lemma yields
\set\begin{equation}\label{eq_up-to-W}
\cJ_WA=I.
\end{equation}
Furthermore, recall that $\cI:=\Ker\,u_A$ contains a
distinguished element $\underline\alpha$ of $W(\bE)$
(lemma \ref{lem_was-third-cond}(ii)); on the other
hand, \eqref{eq_up-to-W} and \eqref{eq_p-lies-deep}
imply that $p\in\cJ_W^tA$, in which case we may apply
lemma \ref{lem_before-name}(i) to get an element
$\underline\alpha':=(\alpha'_n~|~n\in\N)\in\cI$ such
that $\alpha'_0\in\cJ^t$ and $\alpha'_1\in\bE^\times$.
Especially, $\underline\alpha'$ is a distinguished
element of $W(\bE)$, and so the ring
$$
B:=W(\bE)/\underline\alpha' W(\bE)
$$
is perfectoid, for the quotient topology arising from
the natural projection $W(\bE)\to B$ (example
\ref{ex_perfectoid}(ii)). By construction, $u_A$
factors through a surjective and continuous ring
homomorphism
$$
v:B\to A.
$$
Moreover, set $I_B:=\cJ_WB\subset B$, and notice that on
the one hand, $u_B\circ\tau_\bE(\alpha'_0)=pt$ for some
$t\in B^\times$ (lemma \ref{lem_was-third-cond}(iii)),
and on the other hand we have
$u_B\circ\tau_\bE(\alpha'_0)\in I_B^t$, so
\set\begin{equation}\label{eq_p-lies-deep-in-B}
p\in I_B^t.
\end{equation}
Combining with proposition \ref{prop_morel}(ii), we
deduce that $I_B$ is an ideal of definition of $B$,
and taking into account \eqref{eq_up-to-W}, we also get
\set\begin{equation}\label{eq_onto-I}
v(I_B)=I.
\end{equation}
Set $J_B:=I_B^{(p)}$; from \eqref{eq_p-lies-deep-in-B},
\eqref{eq_onto-I} and propositions
\ref{prop_begin-criterion-perf}(ii) and 
\ref{prop_change-topol}(i) we then get a commutative
diagram of graded ring homomorphisms
\set\begin{equation}\label{eq_from-B-to-A}
{\diagram
\gr^\bullet_{I_B}B \ar[r] \ar[d] & \gr^\bullet_IA \ar[d] \\
\gr^\bullet_{J_B}B \ar[r] & \gr^\bullet_JA
\enddiagram}
\end{equation}
both of whose vertical arrows are isomorphisms, and
whose horizontal arrows are induced by $v$.
Since $B$ is complete and separated for its $I_B$-adic
topology, and $A$ is complete and separated for its
$I$-adic topology, the theorem will follow from :

\begin{claim}\label{cl_limit-of-v_n}
The map $v$ induces a ring isomorphism
$v_n:B/I_B^{p^n}\isom A/I^{p^n}$ for every $n\in\N$. 
\end{claim}
\begin{pfclaim}[] We argue by induction on $n$.
For $n=0$, notice that the projection $W(\bE)\to B/I_B$
factors through $\bE$, and by construction, the
kernel of the resulting surjection is $\cJ$. Then
the assertion follows from \eqref{eq_modulo-cJ}.
Next, suppose that $n\geq 0$, and that $v_n$ is an
isomorphism. In this case, $v$ induces isomorphisms
$\gr^r_{I_B}B\isom\gr^r_IA$ for every $r<p^n$.
Due to \eqref{eq_from-B-to-A}, we then deduce
that $v$ also induces isomorphisms
$\gr^r_{J_B}B\isom\gr^r_JA$ for every $r<p^n$, from
which it follows easily that $v$ induces an isomorphism
$B/J^{p^n}_B\isom A/J^{p^n}$. However, clearly
$J\subset I^p$ and $J_B\subset I^p_B$, hence
$J^{p^n}_B\subset I^{p^{n+1}}_B$ and $J^{p^n}\subset I^{p^{n+1}}$,
so finally $v_{n+1}$ is an isomorphism.
\end{pfclaim}
\end{proof}

As a corollary of theorem \ref{th_criterium-perfect}
we obtain the following criterion :

\begin{corollary}\label{cor_BDS}
Let $A$ be a ring, $t\in\N$ any integer, $\beta:=(b_n~|~n\in\N)$
an element of\/ $\bE(A)$; suppose that $A$ is complete and
separated for its $b_0$-adic topology $\cT_0$, and $p=ab_0$
for some $a\in A$. Let also $J:=\bigcup_{n\in\N}b_nA$. Consider
the following conditions :
\begin{enumerate}
\alphaenu
\item
$A/J$ is a perfect $\F_p$-algebra.
\item
The Frobenius endomorphism of $A/b^2_tA$ induces a bijection
$b_{t+1}A/b^2_{t+1}A\isom b_tA/b_t^2A$.
\item
$\Ann_A(b_0)=\Ann_A(b_n)$ for every $n\in\N$.
\item
The Frobenius endomorphism $\Phi_{A/b_iA}$ of $A/b_iA$ induces
an isomorphism
$$
\bar\Phi_{A/b_iA}:A/b_{i+1}A\isom A/b_iA
\qquad
\text{for every $i\in\N$}.
$$
\item
$(A,\cT_0)$ is a perfectoid ring.
\end{enumerate}
Then {\em(e)} is equivalent to the conjunction of\/ {\em(a),(b)}
and {\em(c)}, and also to the conjunction of\/ {\em(c)} and {\em(d)}.
Moreover, if\/ {\em(e)} holds, then $A/\Ann_A(b_0)$ is perfectoid
for the $b_0$-adic topology, and the map $\bar u_A:\bE(A)\to A$
induces an isomorphism of\/ $\bE(A)$-modules :
\set\begin{equation}\label{eq_identify-Anns}
\Ann_{\bE(A)}(\beta)\isom\Ann_A(b_0).
\end{equation}
\end{corollary}
\begin{proof} Suppose first that (a),(b),(c) hold. For every $i,j>0$,
consider the commutative diagram
$$
\xymatrix{
b^{j-1}_iA/b^j_iA \ar[rr]^-{f_{i,j}} \ar[d]_{\phi_{i,j-1}} & &
b^j_iA/b^{j+1}_iA \ar[d]^{\phi_{i,j}} \\
b^{j-1}_{i-1}A/b^j_{i-1}A \ar[rr]^-{f_{i-1,j}} & &
b^j_{i-1}A/b^{j+1}_{i-1}A 
}$$
where $\phi_{i,j}$ is induced by the Frobenius endomorphism
of $A/pA$ for every $i,j\in\N$ with $i>0$, and $f_{i,j}$ is
induced by scalar multiplication by $b_i$, for every
$i,j\in\N$ with $j>0$. We notice :

\begin{claim}\label{cl_two-inductions}
(i)\ \ $f_{i,j}$ is an isomorphism for every $i,j\in\N$
with $j>1$.
\begin{enumerate}
\addenu
\item
$\phi_{i,1}$ is an isomorphism for every $i>t$.
\end{enumerate}
\end{claim}
\begin{pfclaim}(i): The surjectivity of $f_{i,j}$ is clear.
For the injectivity, suppose that $x=b_i^{j-1}y$  and
$b_ix=b_i^{j+1}z$ for some $y,z\in A$. It follows that
$b_i^jy=b_i^{j+2}z$, so $b_i^j(y-b_i^2z)=0$, and therefore
also $b_i(y-b_i^2z)=0$, by (c). Since $j>1$,
we deduce that $x=b_i^jz$, whence the assertion.

(ii): We argue by induction on $i$. The case $i=t+1$ is
(b). Suppose then that $i>t$ and that the surjectivity of
$\phi_{i,1}$ is already known. This means that for every
$x\in A$ there exist $y,z\in A$ such that
$b_{i-1}x=b_{i-1}y^p+b_{i-1}^2z$, {\em i.e.} that
$b_{i-1}(x-y^p-b_{i-1}z)=0$. By assumption (c), we deduce
that $b_i(x-y^p-b_{i-1}z)=0$. Set $z':=b_i^{p-1}z$; then
$b_i(x-y^p-b_iz')=0$, {\em i.e.} $b_ix=b_iy_p+b_i^2z'$,
whence the surjectivity of $\phi_{i+1,1}$.

The injectivity of $\phi_{i,1}$ means the following.
For every $x\in A$ such that $b_{i-1}x^p\in b_{i-1}^2A$
we have $b_ix\in b_i^2A$. Now, let $a\in A$, and suppose
that $b_ia^p\in b_i^2A$; then
$b_{i-1}(b_i^{p-1}a^p)=b_i^{2p-1}a^p\in b_{i-1}^2A$,
so that $b_{i+1}^{2p-1}a=b_ib_{i+1}^{p-1}a\in b_i^2A$.
Say that $b_{i+1}^{2p-1}a=b_i^2y$; then
$$
b_{i+2}^{2p^2-p-1}(b_{i+2}a-b_{i+2}^{p+1}y)=0.
$$
By assumption (c), it follows that
$b_{i+2}^{p-1}(b_{i+2}a-b_{i+2}^{p+1}y)=0$ as well, and
therefore $b_{i+1}a\in b_{i+1}^2A$, hence $\phi_{i+1,1}$
is injective.
\end{pfclaim}

From claim \ref{cl_two-inductions} it follows immediately
that $\phi_{i,j}$ is an isomorphism for every $i>t$ and
every $j>0$. Notice now that $(A,\cT_0)$ is a P-ring, and
$I:=b_iA$ is an ideal of definition of $A$, for every $i>0$,
due to assumptions (a) and (b); moreover, any such $I$
fulfills the conditions of \eqref{subsec_setup-critperf}.
In light of theorem \ref{th_criterium-perfect}, it then
suffices to show that $\bar\Phi_{A/b_iA}$ is bijective for
every $i\geq t$. To this aim, from the snake lemma and
assumption (a), we are easily reduced to checking that
$\Phi_{A/b_{i-1}A}$ induces a bijection $J/b_iA\isom J/b_{i-1}A$
for every $i>t$. Then, since $J$ is the union of its subideals
$b_kA$ (for all $k\in\N$), we are further reduced to checking
that $\Phi_{A/b_{i-1}A}$ induces a bijection
$b_kA/b_iA\isom b_{k-1}A/b_{i-1}A$ for every $k\geq i$.
By considering the filtration
$$
0\subset b_{i+1}A/b_iA\subset b_{i+2}A/b_i\subset\cdots
\subset b_kA/b_iA
$$
an easy induction then reduces to the case where $k=i+1$.
Lastly, a similar induction argument, using the filtration
$0\subset b^{p-1}_{i+1}/b_i\subset b_{i+1}^{p-2}/b_i
\subset\cdots\subset b_{i+1}/b_i$ further reduces the
assertion to the surjectivity of $\phi_{i,j}$ for every
$j=1,\dots,p-1$, which is already known, by the foregoing.

Conversely, if (e) holds, then we get both (b) and (d), by
virtue of proposition \ref{prop_begin-criterion-perf}; then,
notice that the direct limit of the system of maps
$(\bar\Phi_{A/b_iA}~|~i\geq t)$ is the Frobenius endomorphism
of $A/J$, whence (a). Next, let $x\in A$ with $b_0x=0$; it
follows that $(b_nx)^{p^n}=0$ for every $n\in\N$, and since
$A$ is reduced (corollary \ref{cor_perf-are-reduced}(i)),
we get $b_nx=0$, whence (c).

Lastly, suppose that (c) and (d) hold; then we get a
commutative diagram whose horizontal arrows are restrictions
of $\Phi_{A/b_0A}$, and are therefore bijections :
$$
\xymatrix{ b^2_{t+1}A/b_1A \ar[r]^-\sim \ar[d]_{j_{t+1}} &
b^2_tA/b_0A \ar[d]^{j_t} \\
b_{t+1}A/b_1A \ar[r]^-\sim & b_tA/b_0A
}$$
and whose vertical arrows are the natural injections; hence,
we get an induced isomorphism $\Coker\,j_{t+1}\isom\Coker\,j_t$,
whence (b). It has already been remarked that (d) implies (a),
so the proof of the stated equivalences is complete.

Next, suppose that $A$ is perfectoid, so all conditions
(a),$\dots$,(e) hold, and let us notice that
\set\begin{equation}\label{eq_empty-inter}
\Ann_A(b_0)\cap b_nA=0
\qquad
\text{for every $n\in\N$}.
\end{equation}
Indeed, if $x\in A$ and we have $b_0\cdot(b_nx)=0$, then
$b_nx=0$ by (c). This means that $\cT_0$ induces the discrete
topology on $\Ann_A(b_0)$, so $\bar A:=A/\Ann_A(b_0)$ is complete
and separated for its $b_0$-adic topology. Consider now the
commutative diagram :
$$
\xymatrix{ A/b_{t+1}A \ar[r]^-{\bar\Phi_{A/b_iA}} \ar[d]_{f_{1,t+1}} &
A/b_tA \ar[d]^{f_{1,t}} \\
b_{t+1}A/b^2_{t+1}A \ar[r] & b_tA/b^2_tA
}$$
whose left (resp. right) vertical arrow is induced by multiplication
by $b_{t+1}$ (resp. by $b_t$), and whose top (resp. bottom) horizontal
arrow is induced by the Frobenius endomorphism of $A/b_tA$ (resp. of
$A/b^2_tA$). Due to (b),(d) and \eqref{eq_empty-inter}, we see that
$\bar\Phi_{A/b_tA}$ restricts to a bijection
$$
\Ann_A(b_{t+1})=\Ker\,f_{1,t+1}\isom\Ker\,f_{1,t}=\Ann_A(b_t)
$$
and combining with (c) we deduce that the Frobenius endomorphism
of $A/b_iA$ restricts to a bijection $\Ann_A(b_0)\isom\Ann_A(b_0)$
for every $i\in\N$. Then, using (d) we conclude that the Frobenius
endomorphism of $\bar A/b_i\bar A$ induces an isomorphism
$\bar\Phi_{\bar A/b_i\bar A}:\bar A/b_{i+1}\bar A\isom\bar A/b_i\bar A$
for every $i\in\N$. Lastly, it is easily seen that
$\Ann_{\bar A}(b_i)=0$ for every $i\in\N$. Summing up, we have
shown that the ring $\bar A$ fulfills conditions (c) and (d),
so it is perfectoid for its $b_0$-adic topology, by the foregoing.
Now, by virtue of \eqref{eq_empty-inter}, we may naturally identify
$\Ann_A(b)$ with its image in $A/b_0A$, and likewise one easily sees
that $\Ann_{\bE(A)}(\beta)\cap\beta\bE(A)=0$, hence we may identify
$\Ann_{\bE(A)}(\beta)$ with its image in $\bE(A)/\beta\bE(A)$.
By claim \ref{cl_trickier-than-usual}(ii), we have a commutative
diagram of $A_0/b_0A$-modules
$$
\xymatrix{ \bE(A)/\beta\bE(A) \ar[r] \ar[d] &
\beta\bE(A)/\beta^2\bE(A) \ar[d] \\
A/b_0A \ar[r] & b_0A_0/b_0^2A_0
}$$
whose top (resp. bottom) horizontal arrow is induced by
multiplication by $\beta$ (resp. by $b_0$), and whose
vertical arrows are induced by $\bar u_A$. To conclude,
it suffices to remark that the kernel of the top (resp.
bottom) horizontal arrow is $\Ann_{\bE(A)}(\beta)$ (resp.
$\Ann_A(b_0)$).
\end{proof}

Let us point out the following variant of corollary
\ref{cor_BDS} :

\begin{corollary}\label{cor_variant-BDS}
Let $A$ be ring, $b\in A$ an element such that $p\in b^pA$,
and such that the $b$-adic topology $\cT$ of $A$ is complete
and separated. Then $(A,\cT)$ is perfectoid if and only if
the following two conditions hold :
\begin{enumerate}
\alphaenu
\item
The Frobenius endomorphism $\Phi_{A/b^pA}$ of $A/b^pA$ induces
a bijection $A/bA\isom A/b^pA$.
\item
$\Ann_A(b^p)=\Ann_A(b^{p-1})$.
\end{enumerate}
\end{corollary}
\begin{proof} If $A$ is perfectoid, then
$\bar u_{A/pA}:\bE(A)\to A/pA$ is surjective (lemma
\ref{lem_was-third-cond}(i)), so there exists
$(b_n~|~n\in\N)\in\bE(A)$ such that $b_0-b\in pA_0$,
and it follows easily that $b_0/b\in A_0^\times$. Thus,
we may replace $b$ by $b_0$, in which case the assumptions
of corollary \ref{cor_BDS} are fulfilled, and we deduce
that both (a) and (b) hold.

Conversely, if (a) holds, then $A$ is a P-ring, and if
also (b) holds, arguing as in the foregoing we may again
assume that $b=b_0$ for some $(b_n~|~n\in\N)\in\bE(A)$; then
we have, more generally:

\begin{claim}\label{cl_better-estimate}
If (a) holds and there exists $t\in\N$ such that
$\Ann_A(b^p)=\Ann_A(b^p/b_t)$, then $A$ is perfectoid.
\end{claim}
\begin{pfclaim}[]In light of theorem \ref{th_criterium-perfect},
it suffices to show that the Frobenius endomorphism of
$A/b_tA$ induces bijections
$$
\Phi_i:\gr^{(t+1)}_i:=b^i_{t+1}A/b^{i+1}_{t+1}A\isom
\gr^{(t)}_i:=b^i_tA/b^{i+1}_tA
\qquad
\text{for every $i\in\N$}.
$$
However, (a) easily implies that $\Phi_i$ is an isomorphism
for every $i=0,\dots,p^{t+1}-1$. We then argue by induction
on $i\geq p^{t+1}-1$. Thus suppose that $\Phi_i$ is an isomorphism
for such $i$; we consider the commutative diagram
$$
\xymatrix{ \gr_i^{(t+1)} \ar[r]^-{\Phi_i} \ar[d]_{\bar b_{t+1}} &
\gr_i^{(t)} \ar[d]^{\bar b_t} \\
\gr_{i+1}^{(t+1)} \ar[r]^-{\Phi_{i+1}} & \gr_{i+1}^{(t)}
}$$
whose vertical arrows are induced by scalar multiplication
by $b_{t+1}$ and respectively $b_t$. Clearly both vertical
arrows are surjections; thus, we are reduced to checking
that $\bar b_t$ is injective. Hence, say that for some
$a\in A$, the class of $b_t^ia$ in $\gr^{(t)}_i$ lies in the
kernel of $\bar b_t$; this means that there exists $c\in A$
with $b_t^{i+1}a=b_t^{i+2}c$, so that $b_t^{i+1}(a-b_tc)=0$,
and since $i+1\geq p^{t+1}$, our assumption easily implies
that $b_t^i(a-b_tc)=0$ as well. So, the class of $b_t^ia$
vanishes in $\gr^{(t)}_i$, as required.
\end{pfclaim}
\end{proof}

Proposition \ref{prop_change-topol} shows that the
completion of a perfectoid ring for the topology
defined by a $p$-adically open and finitely
generated ideal, is again perfectoid for the $p$-adic
topology. We now wish to prove that the same still
holds for a general finitely generated ideal.

\sset\subsubsection{}\label{subsec_general-complete}
Namely, let $A$ be a perfectoid ring, $I\subset A$
(resp. $\cI\subset\bE:=\bE(A)$) a finitely generated
ideal such that the topology of $A$ (resp. of $\bE$)
agrees with the $I$-adic topology $\cT_I$ (resp. with
the $\cI$-adic topology $\cT_\cI$).
Let also $\beta_\bullet:=(\beta_1,\dots,\beta_n)$ be any
finite system of elements of $\bE$ and $\cJ\subset\bE$
(resp. $J\subset A$) the ideal generated by $\beta_\bullet$
(resp. by $\bar u_A(\beta_1),\dots,\bar u_A(\beta_n)$).
Denote by $\cT_J$ (resp. $\cT_\cJ$) the $J$-adic on $A$
(resp. the $\cJ$-adic topology on $\bE$), and let
$$
(A^\wedge_J,\cT^\wedge_J)
\qquad\text{and}\qquad
(\bE^\wedge_{\!\!\cJ},\cT^\wedge_{\!\!\cJ})
$$
respectively the completion of $(A,\cT_J)$ and of
$(\bE,\cT_\cJ)$. Lastly, let $\cT^\wedge_I$ (resp.
$\cT^\wedge_\cI$) be the $IA^\wedge_J$-adic topology on
$A^\wedge_J$ (resp. the $\cI\bE^\wedge_{\!\!\cJ}$-adic topology
on $\bE^\wedge_{\!\!\cJ}$).

\begin{theorem}\label{th_general-complete}
In the situation of \eqref{subsec_general-complete},
the following holds :
\begin{enumerate}
\item
$(A^\wedge_J,\cT^\wedge_I)$ is perfectoid.
\item
There is a natural isomorphism of topological rings
$\bE(A^\wedge_J,\cT^\wedge_I)\isom
(\bE^\wedge_{\!\!\cJ},\cT^\wedge_\cI)$.
\end{enumerate}
\end{theorem}
\begin{proof} Let $J_\bA\subset\bA:=W(\bE)$ be the ideal
generated by $\tau_\bE(\beta_1),\dots,\tau_\bE(\beta_n)$,
set $J'_\bA:=J_\bA+p\bA$ and denote by $\bA^{\!\wedge}_J$
(resp. $\bA^{\!\wedge}_{J'}$) the $J_\bA$-adic (resp.
$J'_\bA$-adic) completion of $\bA$. Moreover, pick
any finite system of generators $(\gamma_1,\dots,\gamma_k)$
for $\cI$, let $I_\bA\subset\bA$ be the ideal
generated by the system
$(p,\tau_\bE(\gamma_1),\dots,\tau_\bE(\gamma_k))$, and
denote by $\cT^\wedge_{I,\bA}$ indifferently the $I_\bA$-adic
topologies on $\bA^{\!\wedge}_{J'}$ and on $\bA^{\!\wedge}_J$,
and by $\cT^\wedge_{J',\bA}$ the $J'_\bA$-adic topology on
$\bA^{\!\wedge}_{J'}$. With this notation, we have natural
isomorphisms of topological $\bA$-algebras
\set\begin{equation}\label{eq_slowly}
W(\bE^\wedge_{\!\!\cJ},\cT^\wedge_{\!\!\cJ})\isom
(\bA^{\!\wedge}_{J'},\cT^\wedge_{J',\bA})
\qquad\text{and}\qquad
W(\bE^\wedge_{\!\!\cJ},\cT^\wedge_\cI)\isom
(\bA^{\!\wedge}_{J'},\cT^\wedge_{I,\bA})
\end{equation}
(proposition \ref{prop_morel}(ii) and lemma
\ref{lem_Witt-limit}(iv)). On the other hand, since the
$J_\bA$-adic topology is finer than the $J'_\bA$-adic topology,
the completion map $\bA\to\bA^{\!\wedge}_{J'}$ factors through
a unique map of $\bA$-algebras
\set\begin{equation}\label{eq_coarser}
\bA^{\!\wedge}_J\to\bA^{\!\wedge}_{J'}.
\end{equation}

\begin{claim}\label{cl_coarser}
The map \eqref{eq_coarser} is a ring isomorphism.
\end{claim}
\begin{pfclaim} Quite generally, consider any abelian group
$M$, endowed with two topologies $\cT$ and $\cT'$ defined by
cofiltered systems of subgroups $(M_\lambda~|~\lambda\in\Lambda)$
and respectively $(M'_\mu~|~\mu\in\Lambda')$. Endow each
quotient $Q_\lambda:=M/M_\lambda$ with the quotient topology
$\cT'_\lambda$ induced by $\cT'$ via the projection map
$M\to Q_\lambda$, and suppose moreover that $\cT$ is finer
than $\cT'$ and $(Q_\lambda,\cT'_\lambda)$ is separated and
complete for every $\lambda\in\Lambda$. Then it follows
that the natural map
$$
(M,\cT)^\wedge\to(M,\cT')^\wedge
$$
(that extends the completion map $M\to(M,\cT')^\wedge$), is
an isomorphism of abelian groups. Indeed, clearly the family
$(M_\lambda+M'_\mu~|~\lambda\in\Lambda,\mu\in\Lambda')$ is
still a fundamental system of open neighborhoods of zero
in $M$ for the topology $\cT'$, so we have natural isomorphisms
$$
(M,\cT')^\wedge\isom\lim_{(\lambda,\mu)\in\Lambda\times\Lambda'}\,
M/(M_\lambda+M'_\mu)\isom\lim_{\lambda\in\Lambda}\,\lim_{\mu\in\Lambda'}\,
M/(M_\lambda+M'_\mu)\isom\lim_{\lambda\in\Lambda}\,(Q_\lambda,\cT')^\wedge
$$
(where the second isomorphism follows from example
\ref{ex_lim_interchange}(ii)) but
$(Q_\lambda,\cT')^\wedge=Q_\lambda$ under the current
assumptions, whence the assertion.

Recall now that the family
$(W(\cJ^{\La\lambda\Ra})~|~\lambda\in\N[1/p])$ is a fundamental
system of open neighborhoods of zero in for the $J_\bA$-adic
topology on $\bA$ (proposition \ref{prop_morel}(i) and
lemma \ref{lem_mon-fract-powers}(iv)). In view of the
foregoing, the claim will then follow, once we have shown
that the quotients $\bA/W(\cJ^{\La\lambda\Ra})$ are complete
and separated for the $J'_\bA$-adic topology. Moreover,
it is clear that $J'_\bA$-adic topology on
$\bA/W(\cJ^{\La\lambda\Ra})$ agrees with the $p$-adic topology,
so -- taking into account proposition
\ref{prop_Witt-is-complete}(ii) -- we come down to checking
that $W(\cJ^{\La\lambda\Ra})$ is a closed ideal for the $p$-adic
topology $\cT_{p,\bA}$ on $\bA$. But notice that if $\cT_d$
denotes the discrete topology on $\bE$, then
$W(\bE,\cT_d)=(\bA,\cT_{p,\bA})$, so the assertion
is a special case of remark \ref{rem_Witt-limit}(iv).
\end{pfclaim}

\begin{claim} $(\bE^\wedge_{\!\!\cJ},\cT^\wedge_\cI)$ is
perfectoid.
\end{claim}
\begin{pfclaim} Since $\bE$ is a perfect topological ring,
the same holds for $(\bE,\cT_{\!\!\cJ})$, hence also for
$(\bE^\wedge_{\!\!\cJ},\cT^\wedge_{\!\!\cJ})$ (remark
\ref{rem_topology-of-E}(v) and example
\ref{ex_discrete-Witt}(ii)), and by the same token,
$(\bE^\wedge_{\!\!\cJ},\cT^\wedge_\cI)$ is a perfect topological
$\F_p$-algebra. By example \ref{ex_perfectoid}(i), it
then remains only to check that the topology $\cT^\wedge_\cI$
is complete and separated. However, since $\bE$ is complete
and separated for the $\cI$-adic topology, $\bE^\wedge_{\!\!\cJ}$
is complete and separated for the $(\cI+\cJ)$-adic topology
(lemma \ref{lem_double-completion}(ii.b)) and then the
assertion follows from lemma \ref{lem_fontaine}.
\end{pfclaim}

Let $\underline\alpha:=(\alpha_n~|~n\in\N)\in\Ker\,u_A$
be any distinguished element; clearly $\underline\alpha$
is still distinguished in $W(\bE^\wedge_{\!\!\cJ},\cT^\wedge_\cI)$,
hence $W(\bE^\wedge_{\!\!\cJ},\cT^\wedge_\cI)\otimes_\bA A$ is
perfectoid (example \ref{ex_perfectoid}(ii)). From
\eqref{eq_slowly} and claim \ref{cl_coarser} we deduce
that both (i) and (ii) will follow, once we have shown

\begin{claim} There is a natural isomorphism of
topological rings :
$$
(\bA^{\!\wedge}_J,\cT^\wedge_{I,\bA})\otimes_\bA A\isom
(A^{\!\wedge}_J,\cT^\wedge_I).
$$
\end{claim}
\begin{pfclaim}[] By construction, the $J_\bA$-adic topology
induces the $J$-adic topology on $A$, via the projection
$u_A:\bA\to A$. Therefore, $u_A$ extends to a surjective map
$$
u_A^\wedge:\bA^{\!\wedge}_J\to A^{\!\wedge}_J
$$
whose kernel is the topological closure of
$\underline\alpha\bA$ in the $J_\bA$-adic topology of
$\bA^{\!\wedge}_J$ (proposition
\ref{prop_replaces-Mat-Th.8.1}(i,v) and lemma
\ref{lem_still-c-adic}(iv)). On the other hand,
by proposition \ref{prop_morel}(ii) and lemma
\ref{lem_was-third-cond}(i), the $I_\bA$-adic topology
on $\bA$ induces the $I$-adic topology on $A$ via $u_A$,
so the topology $\cT^\wedge_{I,\bA}$ induces the topology
$\cT^\wedge_I$ on $A^\wedge_J$, via $u^\wedge_A$. Thus, we are
reduced to checking that $\underline\alpha\bA^{\!\wedge}_J$
is a closed ideal in the $J_\bA$-adic topology of
$\bA^{\!\wedge}_J$.
By claim \ref{cl_coarser}, it then further suffices to
show that $\underline\alpha\bA^\wedge_{J'}$ is closed for
the topology $\cT^\wedge_{J',\bA}$. But in light
of \eqref{eq_slowly}, the latter assertion follows from
corollary \ref{cor_square-powers}.
\end{pfclaim}
\end{proof}

As an application we get the following generalization
of theorem \ref{th_double-completion}(ii) :

\begin{corollary} Let $(A,\cT)$ be a perfectoid ring,
$\beta_\bullet:=(\beta_1,\dots,\beta_k)$ and
$\gamma_\bullet:=(\gamma_1,\dots,\gamma_l)$ two finite
sequences of elements of $\bE:=\bE(A)$, and denote by
$I$ (resp. $J$) the ideal of $A$ generated by the system
$(\bar u_A(\beta_1),\dots,\bar u_A(\beta_k))$ (resp. by
$\bar u_A(\gamma_1),\dots,\bar u_A(\gamma_l)$). Then the
resulting map \eqref{eq_double-completion} is an isomorphism.
\end{corollary}
\begin{proof} Let $A'$ be the $I$-adic completion
of $A$, and for every $k\in\N$, denote by $(J^k)^c$
the topological closure of the ideal $J^kA'$ relative
to the $IA'$-adic topology; arguing as in the proof of
lemma \ref{lem_double-completion}(ii.b), we are easily
reduced to showing that the $JA'$-adic topology on $A'$
agrees with the linear topology defined by the system
of ideals $((J^k)^c~|~k\in\N)$. 
However, let $K$ be an ideal of definition of $A$, and 
and endow $A'$ with its $KA'$-adic topology $\cT'$;
by theorem \ref{th_general-complete}, the topological
ring $(A',\cT')$ is also perfectoid. Moreover, for every
$i\leq l$, denote by $\gamma'_i\in\bE':=\bE(A')$ the
image of $\gamma_i$ under the homomorphism $\bE\to\bE'$
induced by the completion map $A\to A'$; then clearly
$JA'$ is the ideal generated by the system
$(\bar u_{A'}(\gamma'_1),\dots,\bar u_{A'}(\gamma'_l))$.
Hence, we may replace $A$ by $A'$, after which we may
assume that $A$ is complete and separated for its
$I$-adic topology. To conclude, it then suffices to show :

\begin{claim} In the situation of the corollary, suppose
moreover that $A$ is complete and separated for the
$I$-adic topology. Then the $J$-adic topology of $A$
agrees with the linear topology induced by the system
of ideals $((J^k)^c~|~k\in\N)$.
\end{claim}
\begin{pfclaim}[] From lemma \ref{lem_mon-fract-powers}(iv)
we know already that the $J$-adic topology on $A$ agrees
with the linear topology defined by the system of ideals
$([\gamma_\bullet]^{\La\lambda\Ra}A~|~\lambda\in\N[1/p])$
so it suffices to check that
$$
([\gamma_\bullet]^{\La\lambda\Ra}A)^c\subset
[\gamma_\bullet]^{\La\lambda'\Ra}A
\qquad
\text{for every $\lambda,\lambda'$ in $\N[1/p]$ with
$\lambda'<\lambda$}.
$$
Now, any element of
$([\gamma_\bullet]^{\La\lambda\Ra}A)^c$ can be
written as an $I$-adically convergent series
$\sum_{n\in\N}x_n$, where
$x_n\in[\gamma_\bullet]^{\La\lambda\Ra}A\cap
[\beta_\bullet]^{\La n\Ra}A$ for every $n\in\N$
(again by lemma \ref{lem_mon-fract-powers}(iv)).
Fix $\lambda'',\eps\in\N[1/p]$ with $\lambda'<\lambda''<\lambda$,
and $0<\eps<1-\lambda''/\lambda$; by corollary \ref{cor_bingo}(iii)
we get
$$
x_n\in[\gamma_\bullet]^{\La\lambda''\Ra}\cdot
[\beta_\bullet]^{\La n\eps\Ra}A
\qquad
\text{for every $n\in\N$}.
$$
Now we argue as in the proof of claim \ref{cl_B-is-A-cplt} :
by definition, for every $n\in\N$ there exist
\begin{itemize}
\item
a finite set $S_n\subset\N[1/p]^{\oplus l}$ with
$\mu_1+\cdots+\mu_l=\lambda''$ for every
$\mu:=(\mu_1,\dots,\mu_l)\in S_n$
\item
a system $(a_\mu~|~\mu\in S_n)$ of elements of
$[\beta_\bullet]^{\La n\eps\Ra}A$ such that
$$
x_n=\sum_{\mu\in S_n}a_\mu\cdot\bar u_A(\gamma_\bullet^\mu)
\qquad
\text{where $\gamma_\bullet^\mu:=
\gamma_1^{\mu_1}\cdots\gamma_l^{\mu_l}$ for every
$\mu\in S_n$}.
$$
\end{itemize}
Choose $N\in\N$ such that $\lambda''-\lambda'>lp^{-N}$,
and define $\bar\mu$ and $\mu^*$ as in the proof of
proposition \ref{prop_morel}, so that
$\bar\mu\in S:=\{\nu\in p^{-N}\N^{\oplus l}~|~
\lambda''\geq\nu_1+\cdots+\nu_l>\lambda'\}$
for every $n\in\N$ and every $\mu\in S_n$. It follows that
$$
x=\sum_{\nu\in S}
c_\nu\cdot\bar u_A(\gamma_\bullet^\nu)
\qquad\text{where}\qquad
c_\nu:=
\sum_{n\in\N}\sum_{\substack{ \mu\in S_n \\ 
                            \bar\mu=\nu}}a_\mu\cdot
                       \bar u_A(\gamma_\bullet^{\mu^*}).
$$
Clearly
$\bar u_A(\gamma_\bullet^\nu)\in[\gamma_\bullet]^{\La\lambda'\Ra}A$,
and by lemma \ref{lem_mon-fract-powers}(iv) the series $c_\nu$
converges in the $I$-adic topology of $A$ for every
$\nu\in S$; also it is easily seen that $S$ is a
finite set, whence the contention.
\end{pfclaim}
\end{proof}

\subsection{Homological theory of perfectoid rings}
\label{subsec_Hom-theory-perfectoid}
This section gathers some results on the homological
properties of perfectoid rings. The first result is
a simplification of the criterion of theorem
\ref{th_criterium-perfect}, valid in the case
where a regular sequence generates a special ideal of
definition of the P-ring $A$.

\begin{theorem}\label{th_regular-seq-criterion}
Let $A$ be a P-ring, and $(a_1,\dots,a_k)$ a regular
sequence of elements of $A$ that generates a special
ideal of definition $I$ of $A$. Then the following
conditions are equivalent :
\begin{enumerate}
\alphaenu
\item
The Frobenius endomorphism of $A/pA$ induces an
isomorphism $A/I\isom A/I^{(p)}$.
\item
$A$ is perfectoid.
\end{enumerate}
\end{theorem}
\begin{proof} We know already that (b)$\Rightarrow$(a),
by corollary \ref{cor_jackob}. Suppose that (a) holds.
Then $\bE:=\bE(A)$ is perfectoid (proposition
\ref{prop_third-equiv-cond} and example
\ref{ex_perfectoid}(i)). Arguing as in the proof of
theorem \ref{th_criterium-perfect}, we may then
find elements $\beta_1,\dots,\beta_k$ in $\bE$,
which generate an ideal of definition $\cJ$ of $\bE$,
and such that $\bar u_{A/pA}(\beta_i)$ equals the
image of $a_i$ in $A/pA$, for $i=1,\dots,k$;
moreover, $\bar u_{A/pA}$ induces an isomorphism
\eqref{eq_modulo-cJ}.

Next, since $p\in I^{(p)}$, we may find -- again as
in the proof of theorem \ref{th_criterium-perfect} --
a distinguished element $\underline\alpha'\in\Ker\,u_A$
such that $\alpha'_0\in\Phi_\bE(\cJ)$.
We let $B:=W(\bE)/\underline\alpha' W(\bE)$, and we notice,
as in {\em loc. cit.} that $B$ is perfectoid for the
quotient topology inherited from $W(\bE)$, and $u_A$
factors through a surjective ring homomorphism
$v:B\to A$. We define the ideal $I_B\subset B$ as
in the proof of theorem \ref{th_criterium-perfect},
so \eqref{eq_onto-I} holds as well in the current
situation, and arguing as in {\em loc. cit.} we
also see that $p\in I^{(p)}_B$.
Set $B_0:=B/I_B$ and $A_0:=A/I$; there follows a
commutative diagram
\set\begin{equation}\label{eq_diag-slightly}
{\diagram
\Sym_{B_0}^\bullet(I_B/I^2_B) \ar[r] \ar[d] &
\gr^\bullet_{I_B}B \ar[d] \\
\Sym_{A_0}^\bullet(I/I^2) \ar[r] & \gr^\bullet_IA
\enddiagram}
\end{equation}
where $\gr^\bullet_IA$ is the graded ring associated
to the $I$-adic filtration on $A$, and likewise for
$\gr^\bullet_{I_B}B$. The top horizontal arrow of this
diagram is surjective, and by assumption, the bottom
horizontal arrow is an isomorphism. Furthermore, we
notice :

\begin{claim}\label{cl_slightly-better}
The maps $\bar u_{A/pA}$, $\bar u_{B/pB}$ and $v$ induce
ring isomorphisms
$$
\bE/\cJ^2\isom B/I^2_B\isom A/I^2.
$$
\end{claim}
\begin{pfclaim} We have a commutative diagram
$$
\xymatrix{ \bE/\cJ \ar[r] \ar[d]_{\bar\Phi_{\bE/\cJ}} &
A/I \ar[d]^{\bar\Phi_{A/I}} \\
\bE/\cJ^{(p)} \ar[r] & A/I^{(p)}
}$$
whose horizontal arrows are induced by $\bar u_{A/pA}$.
Now, the vertical arrows are isomorphisms, and we have
already remarked that the top horizontal arrow is also
an isomorphism. Thus, the same holds for the bottom
horizontal arrow; since $\bar u_{A/pA}(\cJ^2)=I^2/pA$, we
deduce already that $\bar u_{A/pA}$ induces an isomorphism
$\bE/\cJ^2\isom A/I^2$. Similarly, since $B$ is perfectoid,
$\bar u_{B/pB}$ induces an isomorphism
$\bE/\alpha'_0\bE\isom B/pB$; however, $\alpha'_0\in\cJ^2$,
$p\in I^2_B$, and $\bar u_{B/pB}(\cJ^2)=I^2_B/pB$, so
$\bar u_{B/pB}$ induces an isomorphism $\bE/\cJ^2\isom B/I^2_B$.
By construction, it is then clear that $v$ induces an
isomorphism $B/I^2_B\isom A/I^2$.
\end{pfclaim}

From claim \ref{cl_slightly-better} we see that $v$
induces an isomorphism $B_0\isom A_0$, and we get as
well isomorphisms of $\bE$-modules
$$
\cJ/\cJ^2\isom I_B/I_B^2\isom I/I^2.
$$
Thus, the left vertical arrow of \eqref{eq_diag-slightly}
is an isomorphism; therefore the same holds for the top
horizontal arrow, and hence also for the right vertical
arrow. Since $B$ is complete and separated for its
$I_B$-adic topology, and $A$ is complete and separated
for its $I$-adic topology, we conclude that $v$ is an
isomorphism of topological rings, whence (b).
\end{proof}

\sset\subsubsection{}\label{subsec_link-with-a-un-etc}
For every $n,r\in\N$, set
$$
P_r:=\N[1/p]^{\oplus r}
\qquad
R_{r,n}:=\Z[P_r]=\Z[T^{1/p^\infty}_1,\dots,T^{1/p^\infty}_r]
$$
so, actually $R_{r,n}$ is independent of $n$, but we also
consider the ring homomorphism
$$
v_{r,n}:R_{r,n}\to R_{r,0}
\qquad
T^a_i\mapsto T_i^{a/p^n}
\qquad
\text{for every $a\in\N[1/p]$ and $i=1,\dots,r$}
$$
for every $r,n\in\N$. Let also $T_\bullet P_r\subset P_r$ be the
ideal generated by $\N^{\oplus r}$, and set
$$
I_{r,n}:=(T_\bullet P_r)\cdot R_{r,n}
\qquad
I^{\lceil s\rceil}_{r,n}:=(T_\bullet P_r)^{\lceil s\rceil}\cdot R_{r,n}
\qquad
\text{for every $s\in\R_+$}
$$
(notation of remark \ref{rem_Witt-are-f-adic}(i)).
Thus $I^{\lceil s\rceil}_{r,n}$ is an ideal of $R_{r,n}$, for
every $s\in\R_+$ and every $r,n\in\N$. Next, let $A$ be
any ring, $u_0:R_{r,0}\to A$ any ring homomorphism, and
for every $n\in\N$ set $u_n:=u_0\circ v_{r,n}:R_{r,n}\to A$,
denote by $\bff^{(n)}$ the sequence
$(u_n(T_1),\dots,u_n(T_r))$ of elements of $A$; for
$n=0$, we shall also write just $\bff$ instead of
$\bff^{(0)}$ and denote by $J$ the ideal generated
by $\bff$. We may then state :

\begin{lemma}\label{lem_generalissimo}
In the situation of \eqref{subsec_link-with-a-un-etc},
suppose furthermore that
$$
\Tor_i^{R_{r,0}}(R_{r,0}/I^{\lceil s\rceil}_{r,0},A)=0
\qquad
\text{for every $i>0$ and every $s\in\R_+$.}
$$
Then the following holds :
\begin{enumerate}
\item
The ring $A$ satisfies condition
$\mathrm{(a)}^\mathrm{un}_{\bff^{(n)}}$ of \eqref{subsec_badabum}
with step $\leq r$, for every $n\in\N$.
\item
$\Ann_A(J^k)=\Ann_A(I^{\lceil 0\rceil}_{r,0}\cdot A)$ for every
$k>0$.
\item
$\Ann_A(J)\cap I^{\lceil 0\rceil}_{r,0}\cdot A=0$.
\end{enumerate}
\end{lemma}
\begin{proof}(i): In case $r=0$, there is nothing to prove,
so we assume henceforth that $r>0$. For every $n\in\N$,
endow $A$ with the $R_{r,n}$-module structure induced by $u_n$.

\begin{claim}\label{cl_reduce-to-f-zero}
$\Tor_i^{R_{r,n}}(R_{r,n}/I^{\lceil s\rceil}_{r,n},A)=0$
for every $i>0$ and every $s\in\R_+$.
\end{claim}
\begin{pfclaim} The base change $v_{r,n}:R_{r,n}\to R_{r,0}$
yields a spectral sequence
$$
E^2_{ij}:=
\Tor^{R_{r,0}}_i(\Tor_j^{R_{r,n}}(R_{r,n}/I^{\lceil s\rceil}_{r,n},R_{r,0}),A)
\Rightarrow\Tor_i^{R_{r,n}}(R_{r,n}/I^{\lceil s\rceil}_{r,n},A).
$$
But since $v_{r,n}$ is an isomorphism and
$v_{r,n}(I^{\lceil s\rceil}_{r,n})=I^{\lceil s/p^n\rceil}_{r,0}$,
we see that $E^2_{ij}=0$ for every $j>0$, and
$E^2_{i0}=\Tor_i^{R_{r,0}}(R_{r,0}/I^{\lceil s/p^n\rceil}_{r,0},A)$.
Thus, $E^2_{i0}=0$ for $i>0$, whence the claim.
\end{pfclaim}

In view of claim \ref{cl_reduce-to-f-zero}, we may
replace $u_0$ by $u_n$, and reduce to checking condition
$\mathrm{(a)^{un}_\bff}$ for $A$. To this aim, define
$A_r$, $I_r$ and the ring homomorphism $\beta_\bff:A_r\to A$
as in \eqref{subsec_badabum}; endow also $R_{r,0}$ with the
$A_r$-module structure induced by the inclusion map
$A_r\to R_{r,0}$, and consider the change of ring spectral
sequence
$$
E(n)^2_{pq}:=\Tor_p^{R_{r,0}}(\Tor_q^{A_r}(A_r/I_r^n,R_{r,0}),A)
\Rightarrow\Tor_{p+q}^{A_r}(A_r/I_r^n,A)
\qquad
\text{for every $n\in\N$}.
$$
It is easily seen that $R_{r,0}$ is a free $A_r$-module,
hence $E(n)^2_{pq}=0$ whenever $q>0$, and
$E(n)^2_{p0}=R_{r,0}/I_{r,0}^n$. There follows a natural
isomorphism
$$
\Tor^{R_{r,0}}_p(R_{r,0}/I^n_{r,0},A)\isom\Tor^{A_r}_p(A_r/I_r^n,A)
\qquad
\text{for every $p,n\in\N$}.
$$
Now, recall that
$I^{n+r}_{r,0}\subset I_{r,0}^{\lceil n+r-1\rceil}\subset I^n_{r,0}$
for every $n\in\N$ (lemma \ref{lem_mon-fract-powers}(iv));
then our assumption implies immediately that the inverse
system $(\Tor^{A_r}_p(A_r/I_r^n,A)~|~n\in\N)$ is uniformly
essentially zero, with step $\leq r$.

(iii): We have just seen that condition
$\mathrm{(a)}^\mathrm{un}_{\bff^{(n)}}$ holds for every
$n$, with step bounded by $r$. For every $n\in\N$,
let $J_n\subset A$ be the ideal generated by $\bff^{(n)}$;
according to remark \ref{rem_logic-argument}(ii), it
follows that there exists a sequence of integers
$(d(i)~|~i\in\N)$ such that inverse system
$(H_{i+1}(\bff^{(n)},J^k_n)~|~k\in\N)$ is uniformly
essentially zero, with step bounded by $d(i)$, for
every $i,n\in\N$. Especially, for $d:=d(r)$ we get
$$
\Ann_A(J_n)\cap J^d_n=0
\qquad
\text{for every $n\in\N$}.
$$
By a simple induction argument, we deduce that
$\Ann_A(J^{k+1}_n)\cap J^d_n=0$ for every $n,k\in\N$.
But clearly, for every $n\in\N$ there exists $k\in\N$
such that $J^{k+1}_n\subset J$, whence
$$
\Ann_A(J)\cap J^d_n=0
\qquad
\text{for every $n\in\N$}
$$
which yields easily the contention.

(ii) is an easy consequence of (iii).
\end{proof}

\begin{proposition}\label{prop_perf-and-Tors}
Let $f:A\to B$ be any homomorphism of perfect $\F_p$-algebras,
$I\subset A$ any ideal. Then the natural morphism
$$
A/I^{\lceil s\rceil}A\derotimes_AB\to
A/I^{\lceil s\rceil}A\otimes_AB[0]
$$
is an isomorphism in $\sD(A\Mod)$, for every $s\in\R_+$
(notation of remark {\em\ref{rem_Witt-are-f-adic}(i)}).
\end{proposition}
\begin{proof} We reduce easily to the case where $I$
is finitely generated, and $f=\bE^*(f_0)$ for a map
$f_0:A_0\to B_0$ of $\F_p$-algebras of finite type
(where $A_0$ and $B_0$ are endowed with their discrete
topologies : notation of \eqref{subsec_ring-struct-on-E}).
Then we may also assume that $I=I_0A$ for some finitely
generated ideal $I_0\subset A_0$, and it suffices to show :

\begin{claim} In the foregoing situation, for every
$\lambda,\lambda'\in\N[1/p]$ with $\lambda'<\lambda$,
the inclusion $I^{\La\lambda\Ra}A\subset I^{\La\lambda'\Ra}A$
induces the zero map
$$
\Tor^A_i(A/I^{\La\lambda\Ra}A,B)\to\Tor^A_i(A/I^{\La\lambda'\Ra}A,B)
\qquad
\text{for every $i>0$}.
$$
\end{claim}
\begin{pfclaim}[] Under the current assumptions, $A$ is
the colimit of the system of $\F_p$-algebras $(A_n~|~n\in\N)$
such that $A_n:=A_0$ for every $n\in\N$, with transition
maps $\phi_n:A_n\to A_{n+1}$ given by the Frobenius map
$\Phi_{A_0}$, and likewise for $B$. Moreover, say that
$\lambda=ap^k$ and $\lambda'=bp^k$ for some $a,b,k\in\N$,
and for every $n\in\N$ define ideals of $A_n$ by the rules :
$$
I_n:=\left\{\begin{array}{ll} 
           0 & \text{for $n<k$} \\
           I_0^{ap^{n-k}} & \text{for $n\geq k$}
       \end{array}\right.
\qquad\text{and}\qquad
I'_n:=\left\{\begin{array}{ll} 
           0 & \text{for $n<k$} \\
           I_0^{bp^{n-k}} & \text{for $n\geq k$}.
       \end{array}\right.
$$
Then it is easily seen that $\phi_n(I_n)\subset I_{n+1}$,
and likewise for $I'_n$, for every $n\in\N$. Furthermore,
the colimit of the system $(I_n~|~n\in\N)$ (resp.
$(I'_n~|~n\in\N)$) is $I^{\La\lambda\Ra}A$ (resp.
$I^{\La\lambda'\Ra}A$). By claim \ref{cl_colim&Tors}(i),
we are therefore reduced to showing that the map induced
by the inclusion
$I_n\subset I'_n$
$$
\Tor_i^{A_0}(A_0/I_n,B_0)\to\Tor_i^{A_0}(A_0/I'_n,B_0)
$$
vanishes, for every sufficiently large $n\in\N$. Since
$a>b$, the latter assertion follows from claim
\ref{cl_vanish-with-bounds}(ii) (applied with
$M:=A$ and $N:=B$).
\end{pfclaim}
\end{proof}

\sset\subsubsection{}\label{subsec_back-to-toids}
We are going now to apply the general setup of
\eqref{subsec_link-with-a-un-etc} to the case where
$A$ is a perfectoid ring, and
$$
\bff:=(\bar u_A(\beta_1),\dots,\bar u_A(\beta_r))
$$
for a given sequence $\beta_\bullet:=(\beta_1,\dots,\beta_r)$
of elements of $\bE:=\bE(A)$. To ease notation, we shall
simply write $P$ instead of $P_r$ and $R$ instead of
$R_{r,0}$. Also, we shall write $I,I^{\lceil s\rceil}\subset R$
instead of $I_{r,0}$, and respectively $I^{\lceil s\rceil}_{r,0}$,
for every $s\in\R_+$. Likewise, we shall write $u$ instead
of $u_0$; hence $u:R\to A$ is the unique ring homomorphism
such that $u(T_i)=\bar u_A(\beta_i)$ for $i=1,\dots,r$.

\begin{proposition}\label{prop_back-to-toids}
With the notation of \eqref{subsec_back-to-toids},
the natural morphism
$$
R/I^{\lceil s\rceil}\derotimes_RA\to
A/I^{\lceil s\rceil}A[0]
$$
is an isomorphism in $\sD(R\Mod)$, for every $s\in\R_+$.
\end{proposition}
\begin{proof} Notice that $R_0:=R\otimes_\Z\F_p$ is a
perfect $\F_p$-algebra, and
$I^{\lceil s\rceil}R_0=I_0^{\lceil s\rceil}R_0$, where
$I_0\subset R_0$ is the ideal generated by
$(T_1,\dots,T_r)$. In light of proposition
\ref{prop_perf-and-Tors}, we deduce
\set\begin{equation}\label{eq_reduce-to-perf}
\Tor^{R_0}_i(R_0/I^{\lceil s\rceil}R_0,\bE)=0
\qquad
\text{for every $s\in\R_+$ and every $i>0$}.
\end{equation}
On the other hand, we have a standard $2$-spectral sequence
$$
E^2_{pq}:=\Tor^{R_0}_p(\Tor^R_q(R/I^{\lceil s\rceil},R_0),\bE)
\Rightarrow\Tor^R_{p+q}(R/I^{\lceil s\rceil},\bE)
$$
(\cite[Th.5.6.6]{We}). However, by construction $I^{\lceil s\rceil}$
is a direct summand of the free $\Z$-module $R$, so $R/I^{\lceil s\rceil}$
is also a free $\Z$-module, and then a standard calculation
shows that $E^2_{pq}=0$ whenever $q>0$. From \eqref{eq_reduce-to-perf}
we see as well that $E^2_{p0}=0$ for every $p>0$, so we conclude
$$
\Tor^R_i(R/I^{\lceil s\rceil},\bE)=0
\qquad
\text{for every $s\in\R_+$ and every $i>0$}.
$$
Next, an easy induction argument using the short exact
sequences of $W(\bE)$-modules
$$
0\to \bE\simeq p^nW_{n+1}(\bE)\to W_{n+1}(\bE)\to W_n(\bE)\to 0
\qquad
\text{for every $n\in\N$}
$$
shows more generally that
\set\begin{equation}\label{eq_isn}
\Tor^R_i(R/I^{\lceil s\rceil},W_n(\bE))=0
\qquad
\text{for every $s\in\R_+$, every $i>0$ and every $n\in\N$}.
\end{equation}
For every $n\in\N$, set $J_n:=\bp_P^{-n}(T_\bullet P)$,
where $\bp_P$ denotes the $p$-Frobenius automorphism
of $P$ (see definition \ref{def_exact-phi}(ii)), and
consider the projection
$$
\pi^{(a,b)}_n:R/J_n^aR\to R/J_n^bR
\qquad
\text{for every $a,b,n\in\N$ with $a\geq b$}.
$$

\begin{claim}\label{cl_resolve-for-Tor}
With the foregoing notation, we have :
\begin{enumerate}
\item
$\Tor_i^R(\pi^{(a,b)}_n,W_m(\bE))=0$\ \ 
for every $a,b,i,m,n\in\N$ with $a-b\geq r$ and $i>0$.
\item
For every $a,n\in\N$ there exists a complex of free
$R$-modules of finite type $L^\bullet_{a,n}$ with a
quasi-isomorphism $L^\bullet_{a,n}\to R/J^a_n[0]$ of
complexes of $R$-modules.
\item
$\Tor_i^R(\pi^{(a,b)}_n,W(\bE))=0$\ \ 
for every $a,b,i,n\in\N$ with $a-b\geq 2r$ and $i>0$.
\end{enumerate}
\end{claim}
\begin{pfclaim} If $r=0$ there is nothing to prove,
so we may assume that $r>0$.
 
(i): Notice that $J_n$ is generated by a
system of $r$ elements of $P$, for every $n\in\N$; in
view of lemma \ref{lem_mon-fract-powers}(ii.c,iv), we
deduce
$$
J_n^a\subset J_n^{\La a\Ra}\subset
(T_\bullet P)^{\lceil (b+r-1)p^{-n}\rceil}\subset J_n^{\La b+r-1\Ra}
\subset J_n^b
$$
so the assertion follows from \eqref{eq_isn}.

(ii): The ideal $J_nR$ is generated by the system
$(T_1^{1/p^n},\dots,T^{1/p^n}_r)$ of elements of the
subring $S_n:=\Z[T_1^{1/p^n},\dots,T^{1/p^n}_r]$;
since $S_n$ is noetherian, we may find a complex
$L'^\bullet_{a,n}$ of free $S_n$-modules of finite type,
with a quasi-isomorphism
$L'^\bullet_{a,n}\to S_n/J^a_nS_n[0]$, and since the
inclusion map $S_n\to R$ is flat, the complex
$L^\bullet_{a,n}:=L'^\bullet_{a,n}\otimes_{S_n}R$ will do.

(iii): For given $a,n\in\N$, let $L^\bullet_{a,n}$ be as
in (ii); it is easily seen that the natural map
$$
\lim_{m\in\N}L^\bullet_{a,n}\otimes_RW_m(\bE)\to
L^\bullet_{a,n}\otimes_RW(\bE)
$$
is an isomorphism in $\sC(R\Mod)$. To ease notation,
set $\cT^a_{i,n,m}:=\Tor^R_i(R/J^a_nR,W_m(\bE))$ for
every $i,m,n\in\N$; in light of \cite[Th.3.5.8]{We}
we deduce natural short exact sequences
$$
0\to R^1\lim_{m\in\N}\cT^a_{i-1,n,m}\to\Tor^R_i(R/J^a_nR,W(\bE))
\to\lim_{m\in\N}\cT^a_{i,n,m}\to 0
$$
for every $i,a,n\in\N$. Then the assertion follows
from (i).
\end{pfclaim}

Now, for any $s\in\R_+$, let $S_s$ be the set of all
pairs $(a,n)\in\N\times\N$ such that $ap^{-n}>s$;
we define a partial ordering on $S_s$ by the rule :
$(a,n)\geq(b,m)$ if and only if $n\geq m$ and
$ap^m\leq bp^n$. With this notation, notice that
$J_m^b\subset J^a_n$ whenever $(a,n)\geq(b,m)$, and
$I^{\lceil s\rceil}$ is the union of its filtered
system of subideals $(J_n^aR~|~(a,n)\in S_s)$.
Taking into account claim \ref{cl_resolve-for-Tor}(iii)
we deduce that
\set\begin{equation}\label{eq_eliminate-m}
\Tor^R_i(R/I^{\lceil s\rceil},W(\bE))=0
\qquad
\text{for every $s\in\R_+$ and every $i>0$}.
\end{equation}
Next, from remark \ref{rem_distinguished}(ii) we get a
short exact sequence of $R$-modules
$$
0\to W(\bE)\xrightarrow{\ \underline\alpha\ }W(\bE)\to A\to 0
$$
where $\underline\alpha$ is any distinguished element
in $\Ker\,u_A$.
Then, by the induced long exact $\Tor$-sequence,
\eqref{eq_eliminate-m} implies already that
$\Tor^R_i(R/I^{\lceil s\rceil},A)$ vanishes for every
$i>1$. By the same token, we also see that
$\Tor^R_1(R/I^{\lceil s\rceil},A)$ is isomorphic to
the kernel of scalar multiplication by $\underline\alpha$
on the $R$-module $W(\bE)/I^{\lceil s\rceil}W(\bE)$;
to conclude the proof, it then suffices to remark :

\begin{claim}
$\underline\alpha W(\bE)\cap I^{\lceil s\rceil}W(\bE)=
\underline\alpha I^{\lceil s\rceil}W(\bE)$.
\end{claim}
\begin{pfclaim}[] Denote by $\cI\subset\bE$ the ideal
generated by $\beta_\bullet$; from proposition \ref{prop_morel}(i)
we easily deduce that
$I^{\lceil s\rceil}W(\bE)=\bigcup_{s'>s}W(\cI^{\lfloor s'\rfloor})$,
and in view of proposition \ref{prop_square-powers}(ii) we get
$$
\underline\alpha W(\bE)\cap I^{\lceil s\rceil}W(\bE)=
\bigcup_{s'>s}(
\underline\alpha W(\bE)\cap W(\cI^{\lfloor s'\rfloor}))=
\bigcup_{s'>s}\underline\alpha W(\cI^{\lfloor s'\rfloor})=
\underline\alpha I^{\lceil s\rceil}W(\bE)
$$
as stated.
\end{pfclaim}
\end{proof}

\begin{corollary}\label{cor_back-to-toids}
In the situation of \eqref{subsec_back-to-toids}, set
$J:=IA$, and denote by $A^\wedge_J$ the $J$-adic completion
of $A$. Let also $C_\bullet$ be any bounded above complex
of flat $A$-modules such that for every $i\in\Z$ there
exists $c(i)\in\N$ with $J^{c(i)}\cdot H_iC_\bullet=0$.
Then the natural morphism
$$
C_\bullet\to A^\wedge_J\otimes_AC_\bullet
$$
is an isomorphism in $\sD(A\Mod)$.
\end{corollary}
\begin{proof} It is an immediate consequence of
lemma \ref{lem_generalissimo}(i), proposition
\ref{prop_back-to-toids} and corollary
\ref{cor_cond-a-and-compl}.
\end{proof}

\begin{proposition}\label{prop_back-and-forth-sec}
In the situation of \eqref{subsec_back-to-toids}, suppose
that the sequence $\beta_\bullet$ generates an open ideal
of\/ $\bE$. Then the following conditions are equivalent :
\begin{enumerate}
\alphaenu
\item
The sequence $\beta_\bullet$ is completely secant in $\bE$.
\item
The sequence
$b_\bullet:=(\bar u_A(\beta_1),\dots,\bar u_A(\beta_r))$
is completely secant in $A$.
\end{enumerate}
\end{proposition}
\begin{proof} Set as usual $J:=IA$, and let $A^\wedge_J$
be the $J$-adic completion of $A$. In view of corollary
\ref{cor_back-to-toids}, we see that the sequence
$b_\bullet$ is completely secant in $A$ if and only
if its image is completely secant in $A^\wedge_J$.
Moreover, if $\cJ\subset\bE$ denotes the ideal generated
by $\beta_\bullet$, we have a natural ring isomorphism
$\bE(A^\wedge_J)\isom\bE^\wedge_\cJ$ (theorem
\ref{th_general-complete}(ii)), and by the same
token the sequence $\beta_\bullet$ is completely
secant in the perfectoid ring $\bE$ if and only
if the same holds for its image in $\bE^\wedge_\cJ$.
Thus, we may replace $A$ by $A^\wedge_J$, and assume
from start that $A$ is $J$-adically complete and separated.
Next, let $(\alpha_n~|~n\in\N)$ be any distinguished
element in $\Ker_,u_A$; since $\cJ$ is open, we have
$\alpha_0^n\in\cJ$ for some integer $n\in\N$, and then
we may find $k\in\N$ such that $\alpha_0$ is contained
in the ideal generated by the sequence
$\beta'_\bullet:=(\beta_1^{1/p^k},\dots,\beta_r^{1/p^k})$.
In light of corollary \ref{cor_Kosz-cptl-sec}, we may
then replace $\beta_\bullet$ with $\beta'_\bullet$, and
assume that $\alpha_0\in\cJ$. Notice now that $\{\cJ^n\}=J^n$
for every $\in\N$ (remark \ref{rem_beta-taut}(viii)), and
the inclusion $\cJ^{n+1}\subset\cJ^n$ is taut for every
$n\in\N$, so theorem \ref{th_taut-one}(ii) yields a
natural isomorphism
$$
\bar\tau:\bigoplus_{n\in\N}\cJ^n/\cJ^{n+1}\isom
\bigoplus_{n\in\N}J^n/J^{n+1}.
$$
Moreover, a simple inspection shows that $\bar\tau$ is an
isomorphism of $A/pA$-algebras; it follows immediately that
the sequence $\beta_\bullet$ is $\bE$-quasi-regular if and
only if the sequence $b_\bullet$ is $A$-quasi-regular (see
remark \ref{rem_augm-of-free}(i)). To conclude, it suffices
now to invoke proposition \ref{prop_quasi-reg-equals-reg}.
\end{proof}

\sset\subsubsection{}\label{subsec_decompose-torsion}
Let $A$ be a perfectoid ring, set $\bE:=\bE(A)$, and
denote by $\underline\alpha:=(\alpha_n~|~n\in\N)$ any
distinguished element of $\Ker\,u_A$.
Set $I_1:=\Ann_\bE(\alpha_0)$, let $I_2\subset\bE$ be
the radical of the ideal $\alpha_0\bE$, and $I_3:=I_1+I_2$;
let also $\bE_i:=\bE/I_i$ for $i=1,2,3$. We claim that
the resulting diagram of ring homomorphisms
$$
\cE
\quad :\quad
{\diagram \bE \ar[r] \ar[d] & \bE_1 \ar[d] \\
\bE_2 \ar[r] & \bE_3
\enddiagram}
$$
is cartesian. Indeed, it is clear the induced map
$\bE\to\bE_1\times_{\bE_3}\bE_2$ is surjective; for
the injectivity, suppose that $x\in I_1\cap I_2$;
then $x^n\in\alpha_0A$ for some $n\in\N$, and
$\alpha_0x=0$, so that $x^{n+1}=0$, and therefore
$x=0$, since $\bE$ is perfect. Notice also that
the $\alpha_0\bE$-adic topology agrees with the
discrete topology on $\bE_2$ and $\bE_3$, and by
arguing as in remark \ref{rem_perfectoid-quot}(ii)
we see that $I_1$ is a closed ideal for the $\alpha_0$-adic
topology of $\bE$. Moreover, since $\bE$ is perfect,
it is easily seen that $I_1$ is a radical ideal; then
it follows that the same holds for $I_3$, since we have
already observed that the latter equals $I_1\oplus I_2$.
Summing up, we deduce that $\bE_i$, endowed with its
$\alpha_0$-adic topology, is a perfectoid ring for
$i=1,2,3$. Consequently, the same holds for the rings
$$
A_i:=A\otimes_{\W(\bE)}W(\bE_i)
\qquad
\text{$i=1,2,3$}
$$
and $\bE(A_i)=\bE_i$, provided we endow each $A_i$ with its
$p$-adic topology (proposition \ref{prop_bar-W-perfect}(i)).
Furthermore, arguing as in the proof of claim
\ref{cl_still-cart} we easily see that the resulting
diagram of ring homomorphisms :
$$
A(\cE)
\quad : \quad
{\diagram A \ar[r] \ar[d] & A_1 \ar[d] \\
A_2 \ar[r] & A_3
\enddiagram
}$$
is still cartesian.

\begin{remark}\label{rem_decompose-torsion}
In the situation of \eqref{subsec_decompose-torsion},
we notice that :

(i)\ \
The image of $\alpha_0$ is regular on $\bE_1$, hence the
induced isomorphism $\bE\isom\bE_1\times_{\bE_3}\bE_2$
identifies $I_1$ with $0\times_{\bE_3}\bE_2$, {\em i.e.}
with $\Ker\,(\bE_2\to\bE_3)$.

(ii)\ \
Therefore, the image of $\bar u_A(\alpha_0)$ is regular
in $A_1$ (proposition \ref{prop_back-and-forth-sec}),
and consequently, also $p$ is regular in $A_1$.

(iii)\ \
For $i=2,3$, the image of $\underline\alpha$ in $W(\bE_i)$
equals $pu$, for some $u\in W(\bE_i)^\times$, so $u_{A_2}$
and $u_{A_3}$ induce isomorphisms of topological rings
$$
\bE_2\isom A_2
\qquad
\bE_3\isom A_3
$$
that identify the projection $\bE_2\to\bE_1$ in $\cE$
with the map $A_2\to A_3$ in $A(\cE)$.

(iv)\ \
Consequently, the isomorphism $A\isom A_1\times_{A_3}A_2$
induced by $A(\cE)$ identifies $\Ann_A(p)$ with
$0\times_{A_3}A_2$, {\em i.e.} with $\Ker\,(A_2\to A_3)$,
and summing up, we get a natural $W(\bE)$-linear
identification
$$
\Ann_\bE(\alpha_0)\isom\Ann_A(p)
\qquad : \qquad
\beta\mapsto\bar u_A(\beta).
$$

(v)\ \
Let $I$ be an ideal of adic definition of $A$, and
$J\subset A$ another ideal such that $IJ=0$. Since
$A$ is reduced (corollary \ref{cor_perf-are-reduced}(i)),
we deduce that $I^n\cap J=0$ for every $n\in\N$, whence
$J+I^n=J\oplus I^n$ for every such $n$. It follows that
$$
\bigcap_{n\in\N}(J+I^n)=\bigcap_{n\in\N}(J\oplus I^n)=
J\oplus\bigcap_{n\in\N}I^n=J
$$
which means that $J$ is closed in the topology of $A$.
\end{remark}

\sset\subsubsection{}\label{subsec_map-on-spec}
Keep the notation of \eqref{subsec_decompose-torsion};
we set
$$
X_\bE:=\Spec\,\bE
\qquad
X_A:=\Spec\,A
\qquad\text{and}\qquad
X_{\bE_i}:=\Spec\,\bE_i
\qquad
X_{A_i}:=\Spec\,A_i
$$
for $i=1,2,3$, and denote by
$$
i_{\bE_k}:X_{\bE_k}\to X_\bE
\qquad
i_{A_k}:X_{A_k}\to X_A
\qquad
\text{for $k=1,2,3$}
$$
the resulting closed immersions. To ease notation, define
$$
\phi:=(\Spec\,\bar u_A):X_A\to X_\bE
$$
and recall that $\phi$ is a continuous and spectral map
(see \eqref{subsec_new-place}).

\begin{proposition}\label{prop_continuous-map}
With the notation of \eqref{subsec_map-on-spec}, the
following holds :
\begin{enumerate}
\item
The map $\phi$ is surjective.
\item
For every quasi-coherent $\cO_{\!X_A}$-module $\cF$ and
every integer $i>0$ we have
$$
R^i\phi_*\cF=0.
$$
\end{enumerate}
\end{proposition}
\begin{proof}(i): Consider any $\fq\in X_\bE$, set
$$
J:=\sum_{x\in\fq}\bar u_A(x)\cdot A
\qquad
S:=\{\bar u_A(x)~|~x\in\bE\setminus\fq\}
$$
and notice that $S$ is a multiplicative subset of $A$.
Let $\fp\in X_A$ be any prime ideal; it is easily
seen that $\phi(\fp)=\fq$ if and only if $S\cap\fp=\emptyset$
and $J\subset\fp$. Thus, we come down to showing that
$JA_S\neq A_S$. Suppose this fails; then there exist
$y\in\bE\setminus\fq$ and sequences
$x_\bullet:=(x_1,\dots,x_k)$, $a_\bullet:=(a_1,\dots,a_k)$
of elements of $\fq$ and respectively $A$, such that
$$
\bar u_A(y)=\sum_{i=1}^k\bar u_A(x_i)\cdot a_i.
$$
Especially, $[y]^{\La 1\Ra}A\subset[x_\bullet]^{\La 1\Ra}A$.
Set $\cJ:=\sum_{i=1}^k\bE x_i\subset\bE$;
by corollary \ref{cor_bingo}(i) we deduce
$$
y\in\cJ^{\La 1-1/p\Ra}\bE\subset\fq
$$
a contradiction.

(ii): Let $f\in\bE$ any element, and set
$[f]:=\bar u_A(f)$; in light of \eqref{eq_was-directly}
we get
\set\begin{equation}\label{eq_directly}
\phi^{-1}(\Spec\,E_f)=\Spec\,A_{[f]}.
\end{equation}
Now, recall that $R^i\phi_*\cF$ is naturally isomorphic
to the sheaf associated to the presheaf
$$
U\mapsto H^i(\phi^{-1}U,\cF)
\qquad
\text{for every open subset $U\subset X_\bE$}.
$$
Then the assertion follows immediately from
\eqref{eq_directly}.
\end{proof}

\sset\subsubsection{}\label{subsec_phi-flat-map}
In the situation of \eqref{subsec_map-on-spec},
let $\underline\alpha:=(\alpha_n~|~n\in\N)$ be any
distinguished element of $\Ker\,u_A$; notice the commutative
diagram of topological spaces
$$
\xymatrix{ X^0_A:=\Spec\,A/pA \ar[rr]^-{\phi_0} \ar[d]_{\iota_A} & &
X^0_\bE:=\Spec\,\bE/\alpha_0\bE \ar[d]^{\iota_\bE} \\
X_A \ar[rr]^-\phi & & X_\bE
}$$
whose vertical arrows are the inclusion maps, and with
$\phi_0:=\Spec\,\omega$, where
$\omega:\bE/\alpha_0\bE\isom A/pA$ is the ring isomorphism
deduced from $u_A\otimes_\Z\F_p$ as in remark
\ref{rem_nice-topology}(ii). Denote also
$$
\cO^\tt_{\!X_\bE}\subset\cO^\tt_{\!X_\bE}
\qquad\text{and}\qquad
\cO^\tt_{\!X_A}\subset\cO^\tt_{\!X_A}
$$
respectively the subsheaf of $\alpha_0$-torsion sections
and the subsheaf of $p$-torsion sections. Notice that
$\cO^\tt_{\!X_\bE}$ is independent of the choice of
$\underline\alpha$ (remark \ref{rem_distinguished}(i));
more precisely, corollary \ref{cor_perf-are-reduced}(i)
implies that
$$
\begin{aligned}
\cO^\tt_{\!X_\bE}=\,& \underline\Gamma_{X^0_\bE}\cO_{\!X_\bE}
=\Ker\,(\cO_{X_\bE}\to i_{\bE_1*}\cO_{X_{\bE_1}}) \\
\cO^\tt_{\!X_A}=\,& \underline\Gamma_{X^0_A}\cO_{\!X_A}
=\Ker\,(\cO_{X_A}\to i_{A_1*}\cO_{X_{A_1}}).
\end{aligned}
$$

\begin{lemma}\label{lem_phi-flat-map}
In the situation of \eqref{subsec_phi-flat-map}, the
following holds :
\begin{enumerate}
\item
The corresponding isomorphism of sheaves
$\phi_0^\flat:\cO_{X^0_\bE}\isom\phi_{0*}\cO_{X^0_A}$ lifts to
a map of sheaves of monoids
$$
\phi^\flat:\cO_{X_\bE}\to\phi_*\cO_{X_A}.
$$
\item
$\phi^\flat$ restricts to an isomorphism of sheaves of
abelian groups
$$
\cO^\tt_{\!X_\bE}\isom\phi_*\cO^\tt_{\!X_A}.
$$
\item
For every open subset $U\subset X_\bE$ we have
$$
\phi^\flat_U(x+y)=\phi^\flat_U(x)+\phi^\flat_U(y)
\qquad
\text{for every $x\in\cO^\tt_{\!X_\bE}(U)$ and
$y\in\cO_{\!X_\bE}(U)$}.
$$
\item
For every closed subset $W_A\subset X^0_A$, the $\cO_{X_A}$-module
$\underline\Gamma_{W_A}\cO_{X_A}$ is quasi-coherent.
\end{enumerate}
\end{lemma}
\begin{proof}(i): For every $f\in\bE$, set $D(f):=\Spec\,\bE_f$
and $[f]:=\bar u_A(f)$, and let
$$
\phi^\flat_f:=(\bar u_A)_f:\bE_f\to A_{[f]}
$$
(notation of \eqref{subsec_localize-mons}). Now, let
$g\in\bE$ be any other element, and suppose that
$D(g)\subset D(f)$; there follows a commutative diagram
$$
\xymatrix{ \bE \ar[d]_{\bar u_A} \ar[r] &
\bE_f \ar[d]_{\phi^\flat_f} \ar[r] &
\bE_g \ar[d]^{\phi^\flat_g} \\
A \ar[r] & A_{[f]} \ar[r] & A_{[g]}
}$$
whose horizontal arrows are the localization maps. Since
all localizations are epimorphisms in the category of
monoids, we deduce that the diagram commutes for every
such $f,g$. Especially, $\phi^\flat_f$ depends only on
$D(f)$ (and not on $f$). Then the rule $D(f)\mapsto\phi^\flat_f$
extends to a morphism of sheaves of monoids on $X_\bE$
(see \cite[Ch.0, \S3.2.3]{EGAI}), and by construction
the resulting diagram
$$
\xymatrix{ \cO_{X_\bE} \ar[rr]^-{\phi^\flat} \ar[d] & &
\phi_*\cO_{X_A} \ar[d] \\
\iota_{\bE*}\cO_{X^0_\bE} \ar[rr]^-{\iota_{\bE*}(\phi_0^\flat)} & &
\iota_{\bE*}\phi_{0*}\cO_{X^0_A}
}$$
commutes.

(iii): It suffices to check the stated identity
for $U$ of the form $\Spec\,\bE_f$, where $f\in\bE$ is
any element. Then, in light of \eqref{eq_annihilate-alpha_0},
we easily reduce to the case where $x,y\in\Img\,(\bE\to\bE_f)$.
Then, we are further reduced to checking that
$\bar u_A(x+y)=\bar u_A(x)+\bar u_A(y)$
for every $x,y\in\bE$ such that $\alpha_0\cdot x=0$.
The latter follows easily from proposition
\ref{prop_combinatorial} and lemma
\ref{lem_was-third-cond}(iii) : details left to the reader.

(ii): Lemma \ref{lem_was-third-cond}(iii) easily implies
that $\phi^\flat$ send $\cO_{X_\bE}$ into $\phi_*\cO_{X_A}$.
In light of (iii), it then suffices to check that the map
$\phi^\flat_f$ restricts to an isomorphism
$\Ann_{\bE_f}(\alpha_0)\isom\Ann_{A_{[f]}}(p)$ for every
$f\in\bE$. However
\set\begin{equation}\label{eq_annihilate-alpha_0}
\Ann_{\bE_f}(\alpha_0)=\Ann_\bE(\alpha_0)_f
\qquad\text{and}\qquad
\Ann_{A_{[f]}}(p)=\Ann_A(p)_{[f]}
\end{equation}
so the assertion follows from remark
\ref{rem_decompose-torsion}(iv).

(iv): Set $W_\bE:=\phi_0(W_A)$; from (ii) we see that
$\phi^\flat$ induces an isomorphism
$$
\underline\Gamma_{W_A}\cO_{X_A}=
\underline\Gamma_{W_A}\circ\underline\Gamma_{X^0_A}\cO_{X_A}
\isom\underline\Gamma_{W_A}
(\iota_{A*}\phi^{-1}_0\iota^{-1}_\bE\underline\Gamma_{X^0_\bE}\cO_{X_\bE})=
\iota_{A*}\phi^{-1}_0\iota^{-1}_\bE\underline\Gamma_{W_\bE}\cO_{X_\bE}
$$
of $\cO_{X_A}$-modules (lemma \ref{lem_adj-Gamma_Z}(v)).
Since the functor $\iota_{A*}$ sends quasi-coherent
$\cO_{X^0_A}$-modules to quasi-coherent $\cO_{X_A}$-modules
(\cite[Ch.I, Cor.9.2.2]{EGAI}), we are then reduced to
checking that $\iota^{-1}_\bE\underline\Gamma_{W_\bE}\cO_{X_\bE}$
is a quasi-coherent $\cO_{X^0_\bE}$-module. However, it
was already remarked in \eqref{subsec_phi-flat-map} that
$\underline\Gamma_{W_\bE}\cO_{X_\bE}=
\underline\Gamma_{W_\bE}\cO^\tt_{X_\bE}$ is an
$\iota_{\bE*}\cO_{X^0_\bE}$-module, so we are further reduced
to checking that $\underline\Gamma_{W_\bE}\cO^\tt_{X_\bE}$ is
a quasi-coherent $\cO_{X_\bE}$-module.

Thus, let $f\in\bE$ be any element, and
$s\in\underline\Gamma_{W_\bE}\cO^\tt_{X_\bE}(D(f))$
any section; pick also a subset $S\subset\bE$ such that
$U_\bE:=X_\bE\setminus W_\bE=\bigcup_{g\in S}D(g)$. By assumption,
$s_{|D(fg)}=0$ for every $g\in S$, and there exist $n\in\N$
and $t\in\bE$ such that $s=f^{-n}t$; it follows that for
every $g\in S$ there exists an integer $m_g\in\N$ such that
$(fg)^{m_g}t=0$ in $\bE$. Since $\bE$ is perfect, we deduce
that $fgt=0$ for every such $g$, hence $(ft)_{|U_\bE}=0$, so
finally $s\in(\Gamma_{W_\bE}\cO^\tt_{X_\bE})_f$, whence the
contention.
\end{proof}

\begin{remark}\label{rem_phi-is-natural}
(i)\ \
Keep the notation of lemma \ref{lem_phi-flat-map}. Notice
that the rule $A\mapsto(\phi,\phi^\flat)$ is functorial
in $A$. Namely, if $g:A'\to A''$ is any continuous ring
homomorphism of perfectoid rings, and $(\phi',\phi'^\flat)$
$(\phi'',\phi''^\flat)$ the pairs attached to $A'$ and
respectively $A''$ as in lemma \ref{lem_phi-flat-map}, then
we have commutative diagrams of topological spaces and of
sheaves of monoids :
$$
\xymatrix{ X_{A''} \ar[r]^-{\phi''} \ar[d]_{\psi_A} &
X_{\bE''} \ar[d]^{\psi_\bE} &
\cO_{X_{\bE'}} \ar[rr]^-{\psi_E^\flat} \ar[d]_{\phi'^\flat} & &
\psi_{\bE*}\cO_{X_{\bE''}} \ar[d]^{\psi_{\bE*}\phi''^\flat}  \\
X_{\!A'} \ar[r]^-{\phi'} & X_{\bE'} &
\phi'_*\cO_{X_{\!A'}} \ar[rr]^-{\phi'_*\psi^\flat_A} & &
\psi_{\bE*}\phi''_*\cO_{X_{A''}}=\phi'_*\psi_{A*}\cO_{X_{A''}}
}$$
where $X_{A'}:=\Spec\,A'$, $X_{\bE'}:=\Spec\,\bE(A')$
and similarly for $X_{A''}$ and $X_{\bE''}$, and
$\psi_A:=\Spec\,g$, $\psi_\bE:=\Spec\,\bE(g)$.
The verification shall be left to the reader.

(ii)\ \
The short exact sequence of $\cO_{\!X_\bE}$-modules
$$
\cE\quad :\quad 
0\to\cO^\tt_{\!X_\bE}\to\cO_{X_\bE}\to i_{\bE_1*}\cO_{X_{\bE_1}}\to 0
$$
induces a commutative ladder for every open subset
$\Omega\subset X_\bE$ :
$$
\xymatrix@C-6pt{
0 \ar[r] & H^0(\Omega,\cO^\tt_{\!X_\bE}) \ar[r] \ar[d]_{\beta^\tt} &
H^0(\Omega,\cO_{X_\bE}) \ar[r] \ar[d]^\beta &
H^0(\Omega,i_{\bE_1*}\cO_{X_{\bE_1}}) \ar[d]^{\beta_1} \ar[r]
& H^1(\Omega,\cO^\tt_{\!X_\bE}) \ar[d]^\gamma \\
0 \ar[r] & H^0(\phi^{-1}\Omega,\cO^\tt_{\!X_A}) \ar[r] &
H^0(\phi^{-1}\Omega,\cO_{X_A}) \ar[r] &
H^0(\phi^{-1}\Omega,i_{A_1*}\cO_{X_{A_1}}) \ar[r]
& H^1(\phi^{-1}\Omega,\cO^\tt_{\!X_A})
}$$
whose rows are the initial terms of the long exact
cohomology sequence attached to $\cE$, and where
$\beta$ is induced by $\phi^\flat$.
The map $\beta_1$ is induced by $i_{\bE_1*}\phi^\flat_1$,
where $\phi_1^\flat:\cO_{X_{\bE_1}}\to\phi_1^*\cO_{X_{A_1}}$
is the morphism associated with the perfectoid ring
$A_1$ and the induced continuous map $\phi_1:X_{\bE_1}\to X_{A_1}$
as in lemma \ref{lem_phi-flat-map}(i). Lastly, $\beta^\tt$
and $\gamma$ are deduced from the restriction
$\cO^\tt_{\!X_\bE}\to\phi_*\cO^\tt_{\!X_A}$ of $\phi^\flat$,
which is an isomorphism of abelian sheaves (lemma
\ref{lem_phi-flat-map}(ii,iii)); furthermore, these
subsheaves are supported on the closed subsets $X^0_\bE$
and respectively $X^0_A$, and $\phi$ restricts to a
homeomorphism $\phi_0:X^0_\bE\isom X^0_A$, so $\beta^\tt$
and $\gamma$ are isomorphisms of abelian groups.
The commutativity of the left square subdiagram is
clear, and that of the central one follows from (i).
To show the commutativity of the right square subdiagram,
pick any $\bar s\in H^0(\Omega,i_{\bE_1*}\cO_{X_\bE})$, and
find a covering $\Omega=\bigcup_{i\in I}\Omega_i$ consisting
of open subsets, such that $\bar s_{|\Omega_i}$ is the image
of some $s_i\in H^0(\Omega_i,\cO_{X_\bE})$, for every
$i\in I$. Set $\Omega_{ij}:=\Omega_i\cap\Omega_j$ for every
$i,j\in I$, and let $s_{ij}:=s_{i|\Omega_{ij}}-s_{j|\Omega_{ij}}$,
and notice that $s_{ij}\in H^0(\Omega_{ij},\cO^\tt_{\!X_\bE})$
for every such $i,j$. Then the image of $\bar s$ in
$H^1(\Omega,\cO^\tt_{X_\bE})$ is the class of the \v{C}ech
cocycle $(s_{ij}~|~i,j\in I)$, and the latter maps to the
class of the cocycle
$(\phi^\flat_{\Omega_{ij}}(s_{ij})~|~i,j\in I)$ in
$H^1(\phi^{-1}\Omega,\cO^\tt_{\!X_A})$. On the other hand,
set $\bar t:=\beta(\bar s)$; for every $i\in I$, the
restriction of $\bar t$ to $\phi^{-1}\Omega_i$ is the
image of
$\phi^\flat_{\Omega_i}(s_i)\in H^0(\Omega_i,\cO_{\!X_A})$,
hence the image of $\bar t$ in
$H^1(\phi^{-1}\Omega,\cO^\tt_{\!X_A})$ is the class of the
cocycle $(t_{ij}~|~i\in I)$ with
$t_{ij}:=t_{i|\phi^{-1}\Omega_{ij}}-t_{j|\phi^{-1}\Omega_{ij}}$
for every $i,j\in I$. Lemma \ref{lem_phi-flat-map}(iii)
implies that $t_{ij}=\phi^\flat_{\Omega_{ij}}(s_{ij})$ for
every $i,j\in I$, whence the assertion.
\end{remark}

The following result shows that if $A$ is a perfectoid
ring, proposition \ref{prop_top-on-opens-fadic-case} holds
even in case the open subset $U$ of $\Spec\,A$ is not
quasi-compact.

\begin{lemma}\label{lem_topology-on-opens}
Let $A$ be a perfectoid ring, and $U\subset X_A:=\Spec\,A$
an open subset containing the analytic locus of $X_A$.
We have :
\begin{enumerate}
\item
There exists a unique topology $\cT_U$ on
$A_U:=\cO_{X_A}(U)$ such that $(A_U,\cT_U)$ is f-adic
and the restriction map $\rho_U:A\to A_U$ is open.
\item
$(A_U,\cT_U)$ is complete and separated.
\item
Let $f:A\to B$ be any morphism of f-adic topological
rings, such that the image of\/ $\Spec\,f$ lies in $U$.
Then the resulting map $f_U:(A_U,\cT_U)\to B$ is continuous.
\end{enumerate}
\end{lemma}
\begin{proof}(i): We pick a system of elements
$(\beta_1,\dots,\beta_r)$ of $\bE(A)$ such that the
system $(f_i:=\bar u_A(\beta_i)~|~i=1,\dots,r)$ generates
an ideal $I$ of adic definition of $A$. Like in the proof
of proposition \ref{prop_top-on-opens-fadic-case}(i),
there exists a unique topology $\cT_U$ on $A_U$ such that
$(A_U,\cT_U)$ is a topological group for its additive
group structure, and we need to check the following.
For every $a\in A_U$ and every $n\in\N$ there exists
$k\in\N$ such that
$$
a\cdot\rho_U(I^k)\subset\rho_U(I^n).
$$
Now, for $i=1,\dots,r$, let $a_i\in\Spec\,A_{f_i}$ be the
image of $a$ under the restriction map
$\rho_i:A_U\to A_{f_i}$; we may find $m\in\N$ and
$b_1,\dots,b_r\in A$ such that $a_i=f_i^{-m}b_i$, hence
$f_i^ma-b_i\in\Ker\,\rho_i$ for $i=1,\dots,r$. Pick a
family $(g_\lambda~|~\lambda\in\Lambda)$ of elements of
$A$ such that $U=\bigcup_{\lambda\in\Lambda}V_\lambda$, where
$V_\lambda:=\Spec\,A_{g_\lambda}$ for every $\lambda\in\Lambda$.
Notice that $(f_i^ma-b_i)_{|V_\lambda}$ lies in the kernel
of the restriction map
$\cO_U(V_\lambda)\to\cO_U(V_\lambda\cap\Spec\,A_{f_i})$,
for every $i=1,\dots,r$ and every $\lambda\in\Lambda$. Thus,
for every such $\lambda$ and $i$ there exist
$t_\lambda,s_\lambda\in\N$ and $c_{i,\lambda}\in A$ such that
$$
g_\lambda^{s_\lambda}\cdot(f_i^ma-b_i)_{|V_{\lambda}}=
\rho_\lambda(c_{i,\lambda})
\qquad\text{and}\qquad
f_i^{t_\lambda}c_{i,\lambda}=0
\qquad\text{in $A$}
$$
where $\rho_\lambda:A\to A_{g_\lambda}$ is the localization map.
From corollary \ref{cor_perf-are-reduced}(i) we then deduce
that $f_ic_{i,\lambda}=0$, so
$$
(f_i^{m+1}a-f_ib_i)_{|V_\lambda}=0
\qquad
\text{for every $\lambda\in\Lambda$ and every $i=1,\dots,r$}
$$
{\em i.e.} $f_i^{m+1}a=\rho_U(f_ib_i)$ for every $i=1,\dots,r$.
From this, we may argue as in the proof of proposition
\ref{prop_top-on-opens-fadic-case}(i), to derive the
assertion.

(ii): It suffices to show that $\Ker\,\rho_U$ is a closed ideal
of $A$. However, clearly $I\cdot\Ker\,\rho_U=0$, so the assertion
follows from remark \ref{rem_decompose-torsion}(v).

(iii) is proven as proposition
\ref{prop_top-on-opens-fadic-case}(iii).
\end{proof}

We wish next to show that the formula for the Teichm\"uller
representative of a sum obtained in proposition
\ref{prop_combinatorial}, is still valid in f-adic rings
of the type considered in lemma \ref{lem_topology-on-opens}.

\sset\subsubsection{}\label{subsec_new-formula}
Indeed, resume the situation of \eqref{subsec_phi-flat-map},
and let $U\!\subset\!X_\bE$ be an open subset containing
$X_\bE\setminus\Spec\,\bE/\bE^{\circ\circ}$; from
\eqref{eq_directly} we easily see that
$X_A\setminus\Spec\,A/A^{\circ\circ}\subset\phi^{-1}U$, so
$$
\bE_U:=\cO_{X_\bE}(U)
\qquad\text{and}\qquad
A_U:=\cO_{X_A}(\phi^{-1}U)
$$
admit the f-adic topologies $\cT_U$ and respectively
$\cT_{\phi^{-1}U}$ furnished by lemma \ref{lem_topology-on-opens}.

\begin{proposition}\label{prop_new-formula}
With the notation of \eqref{subsec_new-formula}, the
following holds :

{\em (i)}\ \
The map $\phi_U^\flat:(\bE_U,\cT_U)\to(A_U,\cT_{\phi^{-1}U})$
is continuous.

{\em (ii)}\ \
For every $\underline x:=(x_0,\dots,x_k)\in\bE_U^{k+1}$ we
have the identity
$$
\phi^\flat_U(x_0+\cdots+x_k)=
\sum_{n\in\N}p^n\cdot\sum_{\sigma\in\Sigma^{(k)}_n}
c_\sigma\cdot\phi^\flat_U(\underline x^\sigma)
$$
where $\Sigma^{(k)}_n$ is defined as in
\eqref{subsec_combinatorial-Teich} for every $n\in\N$, and
$c_\sigma\in\Z_p$ is provided by proposition
{\em\ref{prop_combinatorial}}, for every
$\sigma\in\bigcup_{n\in\N}\Sigma^{(k)}_n$; the convergence
of the series is relative to the topology $\cT_{\phi^{-1}U}$.
\end{proposition}
\begin{proof} (Recall that
$\underline x^\sigma:=x_0^{\sigma_0}\cdots x_k^{\sigma_k}$
for every $\sigma:=(\sigma_0,\dots,\sigma_k)\in\N[1/p]^{k+1}$).

(ii): Let $(\alpha_n~|~n\in\N)$ be any distinguished element
in $\Ker\,u_A$; since $\alpha_0$ is topologically nilpotent
in $\bE$, its image is topologically nilpotent in $\bE_U$,
hence there exists $c\in\N$ such that
$\alpha_0^c\cdot x_i\in\Img(\bE\to\bE_U)$ for $i=0,\dots,k$.
In light of lemma \ref{lem_was-third-cond}(iii), we deduce
that
$$
y_n(\underline x):=\sum_{\sigma\in\Sigma_n^{(k)}}
p^n\cdot c_\sigma\cdot\phi^\flat_U(x^\sigma)\in
p^{n-c}\cdot\Img(A\to A_U)
\qquad
\text{for every integer $n\geq c$}.
$$
Taking into account lemma \ref{lem_topology-on-opens}(ii),
it follows that the series $\sum_{n\in\N}y_n(\underline x)$
converges for the topology $\cT_{\phi^{-1}U}$ to a well
determined element $y\in A_U$, and it remains only to check
that $y=\phi^\flat_U(x_0+\cdots+x_k)$. To this aim, pick a
family $(g_\lambda~|~\lambda\in\Lambda)$ of elements of $\bE$,
such that $U=\bigcup_{\lambda\in\Lambda}V_\lambda$, where
$V_\lambda:=\Spec\,\bE_{g_\lambda}$ for every $\lambda\in\Lambda$.
Then, for every such $\lambda$ we may find $n_\lambda\in\N$
such that $(g_\lambda^{n_\lambda}\cdot x_i)_{|V_\lambda}\in
\Img(\bE\to\bE_{g_\lambda})$ for $i=0,\dots,k$, and it suffices
to check that
$([g^{n_\lambda}_\lambda]\cdot\phi^\flat_U(x_0+\cdots+x_k))_{|\phi^{-1}V_\lambda}=
([g^{n_\lambda}_\lambda]\cdot y)_{|\phi^{-1}V_\lambda}$, where
$[g^{n_\lambda}_\lambda]:=\bar u_A(g^{n_\lambda}_\lambda)$
for every such $\lambda$. However, we have
$$
[g^{n_\lambda}_\lambda]\cdot\phi^\flat_U(x_0+\cdots+x_k)=
\phi^\flat_U(g^{n_\lambda}_\lambda\cdot x_0+\cdots+
g^{n_\lambda}_\lambda\cdot x_k)
\qquad\text{and}\qquad
[g^{n_\lambda}_\lambda]\cdot y=\sum_{n\in\N}
y_n(g^{n_\lambda}_\lambda\cdot\underline x)
$$
where $g^{n_\lambda}_\lambda\cdot\underline x:=
(g^{n_\lambda}_\lambda\cdot x_0,\dots,g^{n_\lambda}_\lambda\cdot x_k)$.
Hence, we may replace $\underline x$ by
$g^{n_\lambda}_\lambda\cdot\underline x$, after which we may
assume that $x_i\in\Img(\rho_U:\bE\to\bE_U)$ for $i=0,\dots,k$.
In this case, say that $x_i=\rho_U(z_i)$ for $i=0,\dots,k$
and for certain $z_0,\dots,z_k\in\bE$; we have
$$
\phi^\flat_U(x_0+\cdots+x_k)=
\rho_{\phi^{-1}U}\circ\phi^\flat_{X_\bE}(z_0+\cdots+z_k)=
\rho_{\phi^{-1}U}\circ\bar u_A(z_0+\cdots+z_k)
$$
where $\rho_{\phi^{-1}U}:A\to A_U$ is the corresponding
restriction map; on the other hand, since $\rho_{\phi^{-1}U}$
is continuous for the topology $\cT_{\phi^{-1}U}$, we have as
well
$$
\begin{aligned}
y=\, & \sum_{n\in\N}p^n\cdot\sum_{\sigma\in\Sigma^{(k)}_n}
c_\sigma\cdot\rho_{\phi^{-1}U}\circ\phi^\flat_{X_\bE}(\underline z^\sigma)
=\rho_{\phi^{-1}U}\Bigl(\sum_{n\in\N}p^n\cdot\sum_{\sigma\in\Sigma^{(k)}_n}
c_\sigma\cdot\phi^\flat_{X_\bE}(\underline z^\sigma)\Bigr) \\
=\, & \rho_{\phi^{-1}U}\Bigl(\sum_{n\in\N}p^n\cdot\sum_{\sigma\in\Sigma^{(k)}_n}
c_\sigma\cdot\bar u_A(\underline z^\sigma)\Bigr)
=\rho_{\phi^{-1}U}\circ u_A
\Bigl(\sum_{n\in\N}p^n\cdot\sum_{\sigma\in\Sigma^{(k)}_n}
c_\sigma\cdot\tau_A(\underline z^\sigma)\Bigr)
\end{aligned}
$$
so the sought identity follows from proposition
\ref{prop_combinatorial}.

(i): Fix $x_0\in\bE_U$; the proof of (ii) shows that
for every $k\in\N$ there exists $n\in\N$ such that
$y_r(x_0+x)\in p^k\cdot\Img\,(A\to A_U)$ for every
$r\geq n$ and every $x\in\Img\,(\bE\to\bE_U)$. On
the other hand, clearly the mappings $y_i$ are
continuous for every $i\in\N$, relative to the
topologies $\cT_U$ and $\cT_{\phi^{-1}U}$. The assertion
follows immediately.
\end{proof}

\sset\subsubsection{}\label{subsec_large-open-subset}
Keep the notation of \eqref{subsec_phi-flat-map}, let
$U_\bE\subset X_\bE$ be any open subset, and set
$U_A:=\phi^{-1}U_\bE$, $Z_\bE:=X_\bE\setminus U_\bE$,
$Z_A:=X_A\setminus U_A$. Denote by $j_\bE:U_\bE\to X_\bE$
and $j_A:U_A\to X_A$ the open immersions; the map $\phi^\flat$
of lemma \ref{lem_phi-flat-map}(i) induces a morphism of
sheaves of monoids
$$
\psi:j_{\bE*}\cO_{U_\bE}
\xrightarrow{\ j_{\bE*}j_\bE^*(\phi^\flat)\ }
j_{\bE*}j^*_\bE\phi_*\cO_{X_A}\xrightarrow{\ \sim\ }
\phi_*j_{A*}\cO_{U_A}.
$$
Moreover, let
$$
\bar\cO_{X_\bE}:=\Img\,(\cO_{X_\bE}\to j_{\bE*}\cO_{U_\bE})
\qquad
\bar\cO_{X_A}:=\Img\,(\cO_{X_A}\to j_{A*}\cO_{U_A}).
$$
Denote by
$\pi_\bE:j_{\bE*}\cO_{U_\bE}\to R^1\underline\Gamma_{Z_\bE}\cO_{X_\bE}$
and $\pi_A:j_{A*}\cO_{U_A}\to R^1\underline\Gamma_{Z_A}\cO_{X_A}$
the projections, set $q:=\bar u_A(\alpha_0)$, and for
every $c\in\N[1/p]$ define
$$
\begin{aligned}
\alpha_0^{-c}\cO_{X_\bE}:=\,& \Ker\,\alpha_0^c\cdot\pi_\bE
& \qquad & &
p^{-c}\cO_{X_A}:=\,& \Ker\,q^c\cdot\pi_A \\
\cQ^c_{U_\bE}:=\,& \pi_\bE(\alpha_0^{-c}\cO_{X_\bE})
& \qquad & &
\cQ^c_{U_A}:=\,& \pi_A(p^{-c}\cO_{X_A}).
\end{aligned}
$$
Also, set
$$
\cR^c:=\Coker\,(\phi_*\bar\cO_{X_A}\to\phi_*(p^{-c}\cO_{X_A}))
\qquad\text{and let}\qquad
\lambda:\cR^c\to\phi_*\cQ^c_{U_A}
$$
be the natural morphism of $\cO_{X_A}$-modules; lemma
\ref{lem_was-third-cond}(iii) easily implies that
$\psi$ restricts to a morphism of sheaves
of sets
$$
\psi^c:\alpha_0^{-c}\cO_{X_\bE}\to\phi_*(p^{-c}\cO_{X_A})
\qquad
\text{for every $c\in\N[1/p]$}.
$$

\begin{proposition}\label{prop_Q-c-tricky}
With the notation of \eqref{subsec_large-open-subset},
suppose moreover that :
\begin{enumerate}
\alphaenu
\item
$\Spec\,\bE[\alpha_0^{-1}]\subset U_\bE$.
\item
$U_\bE$ is quasi-compact.
\item
$c\leq p/(p-1)$.
\end{enumerate}
Then we have :
\begin{enumerate}
\item
The map $\psi^c$ descends to a morphism of sheaves of sets
$$
\bar\psi{}^c:\cQ^c_{U_\bE}\to\cR^c
\xrightarrow{\ \lambda\ }\phi_*\cQ^c_{U_A}.
$$
\item
There exist unique morphisms of sheaves of sets
$$
\bar\bp^c_\bE:\cQ^{c/p}_{U_\bE}\to\cQ^c_{U_\bE}
\qquad\text{and}\qquad
\bar\bp^c_A:\cQ^{c/p}_{U_A}\to\cQ^c_{U_A}
$$
fitting into commutative diagrams
$$
\xymatrix{
\alpha_0^{-c/p}\cO_{X_\bE} \ar[r] \ar[d] &
\alpha^{-c}_0\cO_{X_\bE} \ar[d] &
p^{-c/p}\cO_{X_A} \ar[r] \ar[d] & p^{-c}\cO_{X_A} \ar[d] \\
\cQ^{c/p}_{U_\bE} \ar[r]^-{\bar\bp^c_\bE} & \cQ^c_{U_\bE} &
\cQ^{c/p}_{U_A} \ar[r]^-{\bar\bp^c_A} & \cQ^c_{U_A}
}$$
whose vertical arrows are the projections, and whose top
horizontal arrows are the restrictions of the $p$-Frobenius
endomorphisms of the sheaves of monoids $j_{\bE*}\cO_{U_\bE}$
and respectively $j_{A*}\cO_{U_A}$ (with composition law
given by multiplication).
\item
Moreover, we have a commutative diagram of sheaves of sets :
$$
\xymatrix{
\cQ^{c/p}_{U_\bE} \ar[rr]^-{\bar\bp^c_\bE} \ar[d]_{\bar\psi{}^{c/p}}
& & \cQ^c_{U_\bE} \ar[d]^{\bar\psi{}^c} \\
\phi_*\cQ^{c/p}_{U_A} \ar[rr]^-{\phi_*(\bar\bp^c_A)} & &
\phi_*\cQ^c_{U_A}.
}$$
\item
If $c\leq 1$, the map $\bar\psi{}^c$ is an isomorphism
of sheaves of abelian groups.
\end{enumerate}
\end{proposition}
\begin{proof}(i): Let $f\in\bE$ be any element, and set
$V:=\Spec\,\bE_f$; we easily come down to showing

\begin{claim}\label{cl_difference}
If $c\leq p/(p-1)$, we have
$$
\psi^c_V(x+y)-\psi^c_V(x)\in\bar\cO_{X_A}(\phi^{-1}V)
\qquad
\text{for every $x\in\alpha_0^{-c}\cO_{X_\bE}(V)$ and
$y\in\bar\cO_{X_\bE}(V)$}.
$$
\end{claim}
\begin{pfclaim} Since $U_\bE$ is quasi-compact,
$j_{\bE*}\cO_{\!X_\bE}$ is a quasi-coherent $\cO_{X_\bE}$-module
(\cite[Ch.I, Cor.9.2.2]{EGAI}), hence the same holds for
$\alpha_0^{-c}\cO_{X_\bE}$, and consequently
$$
\alpha_0^{-c}\cO_{X_\bE}(V)=\alpha_0^{-c}\cO_{X_\bE}(X_\bE)_f.
$$
We may therefore find $n\in\N$ such that $f^nx$
(resp. $f^ny$) is the restriction to $V$ of a global
section $x'\in\alpha_0^{-c}\cO_{X_\bE}(X_\bE)$ (resp.
of an element $y'\in\bE$), and it suffices to
show that
$$
\psi^c_{X_\bE}(x'+y')-\psi^c_{X_\bE}(x')\in
\Img\,(A\to j_{A*}\cO_{U_A}(X_A)=\cO_{X_A}(U_A)).
$$
Hence, we may replace $V$ by $X_\bE$, and assume from
start that
$$
x\in\cO_{X_\bE}(U_\bE)
\qquad\text{and}\qquad
\alpha_0^c\cdot x,y\in\Img\,(\bE\to\cO_{X_\bE}(U_\bE))
$$
and we have to show that
$\phi^\flat_{U_\bE}(x+y)-\phi^\flat_{U_\bE}(x)\in
\Img\,(A\to\cO_{X_A}(U_A))$.
Notice that the assertion is independent of the topology
$\cT$ of $A$ such that $(A,\cT)$ is perfectoid; by virtue
of proposition \ref{prop_change-topol}(i) we may then
assume that $\cT$ is the $p$-adic topology on $A$, and
therefore the topology of $\bE$ agrees with the
$\alpha_0$-adic topology. In this case, assumption (a),
proposition \ref{prop_new-formula}(ii) and lemma
\ref{lem_was-third-cond}(iii) imply that
$\phi^\flat_{U_\bE}(x+y)-\phi^\flat_{U_\bE}(x)-\phi^\flat_{U_\bE}(y)$
is a $p$-adically convergent series $\sum_{n>0}z_n$, where
$z_n$ is a $\Z_p$-linear combination of terms of the form
$\phi^\flat_{U_\bE}(\alpha_0^n\cdot x^\sigma\cdot y^{1-\sigma})$,
with $1-\sigma,\sigma\in p^{-n}\N\setminus p^{1-n}\N$, for
every $n\in\N$. Hence, it suffices to check that
$n\geq p\cdot\sigma/(p-1)$ for every such $\sigma$ and
every integer $n>0$. We leave the easy verification to
the reader.
\end{pfclaim}

(ii): The uniqueness of $\bar\bp^c_{U_A}$ is clear.
For the existence, let $V$ be as in the foregoing;
we come down to showing that
$$
(x+y)^p-x^p\in\bar\cO_{X_A}(V)
\qquad
\text{for every $x\in p^{-c/p}\cO_{X_A}(V)$ and every
$y\in\bar\cO_{X_A}(V)$}.
$$
However,
$(x+y)^p-x^p=y^p+\sum_{i=1}^{p-1}\binom{p}{i}\cdot x^iy^{p-i}$,
and $px^i\in\bar\cO_{X_A}(V)$ for every $i=1,\dots,p-1$,
since $c\leq p/(p-1)$, whence the contention. Assertion
(iii) shall be left to the reader.

(iv): Let $V$ be as in the foregoing; in order to prove
that $\bar\psi{}^c$ is a morphism of abelian sheaves when
$c\leq 1$, we come down to showing that
$$
\psi^1_V(x+y)-\psi^1_V(x)-\psi^1_V(y)\in
\bar\cO_{X_A}(\phi^{-1}V)
\qquad
\text{for every $x,y\in\alpha_0^{-1}\cO_{X_\bE}(V)$}
$$
and arguing as in the proof of claim \ref{cl_difference},
we reduce to the case where $V=X_\bE$ and $x,y$ are
sections of $\cO_{X_\bE}(U_\bE)$ such that
$\alpha_0x,\alpha_0y\in\Img\,(\bE\to\cO_{X_\bE}(U_\bE))$.
Then, it suffices to check that
$$
\phi^\flat_{U_\bE}(\alpha_0^n\cdot y^\sigma\cdot x^{1-\sigma})
\in\Img(A\to\cO_{X_A}(U_A))
\qquad
\text{whenever $n>0$ and $1-\sigma,\sigma\in\N[1/p]$}
$$
which is clear. Next notice that, by remark
\ref{rem_phi-is-natural}, the pair $(\phi,\phi^\flat)$
restricts to isomorphisms of schemes
$$
(\phi_2,\phi_2^\flat):X_{A_2}\isom X_{\bE_2}
\qquad
(\phi_3,\phi^\flat_3):X_{A_3}\isom X_{\bE_3}
\qquad\text{such that}\qquad
\phi_2(Z_A)=Z_\bE.
$$

\begin{claim}\label{cl_sweat-and-tears}
Assertion (iv) holds for $U_\bE=\Spec\,\bE[\alpha_0^{-1}]$.
\end{claim}
\begin{pfclaim} We have commutative diagrams
$$
\xymatrix{ U_{\bE_1}:=\Spec\,\bE_1[\alpha_0^{-1}]
\ar[r] \ar[d]_{j_{\bE_1}} & U_\bE \ar[d]^{j_\bE} &
U_{A_1}:=\Spec\,A_1[p^{-1}]
\ar[r] \ar[d]_{j_{A_1}} & U_A \ar[d]^{j_A} \\
X_{\bE_1} \ar[r]^-{i_{\bE_1}} & X_\bE &
X_{A_1} \ar[r]^-{i_{A_1}} & X_A
}$$
whose top horizontal arrows are isomorphisms of
schemes, whence commutative diagrams of abelian
sheaves
$$
\xymatrix{
\bar\cO_{X_\bE} \ar[r]^-\sim \ar[d] & i_{\bE_1*}\cO_{X_{\bE_1}} \ar[d] &
\bar\cO_{X_A} \ar[r]^-\sim \ar[d] & i_{A_1*}\cO_{X_{A_1}} \ar[d] \\
j_{\bE*}\cO_{U_\bE} \ar[r]^-\sim & i_{\bE_1*}\circ j_{\bE_1*}\cO_{U_{\bE_1}} &
j_{A*}\cO_{U_A} \ar[r]^-\sim & i_{A_1*}\circ j_{A_1*}\cO_{U_{A_1}}
}$$
whose left (resp. right) vertical arrows are the inclusion
maps (resp. are induced by the inclusion maps
$\cO_{X_{\bE_1}}\to j_{\bE_1*}\cO_{U_{\bE_1}}$ and
$\cO_{X_{A_1}}\to j_{A_1*}\cO_{U_{A_1}}$). Moreover, denote by
$\phi_1:X_{A_1}\to X_{\bE_1}$ the continuous map associated
with $A_1$, as in \eqref{subsec_map-on-spec}, and by
$$
\psi_1:j_{\bE_1*}\cO_{U_{\bE_1}}\to\phi_{1*}j_{A_1*}\cO_{U_{A_1}}
\qquad
\text{(resp.\
$\bar\psi{}^c_1:\cQ^c_{U_{\bE_1}}\to\phi_{1*}\cQ^c_{U_{A_1}}$\ )}
$$
the corresponding morphism of sheaves of monoids as in
\eqref{subsec_large-open-subset} (resp. morphism of sheaves
of sets, as in (i)); then $i_{\bE_1}\circ\phi_1=\phi\circ i_{A_1}$
(remark \ref{rem_phi-is-natural}), and the foregoing
isomorphisms identify $\psi$ with $i_{\bE_1*}\psi_1$.
So finally, $\bar\psi{}^c$ is identified with
$i_{\bE_1*}\bar\psi{}^c_1$ for every $c\leq p/(p-1)$,
and we may replace from start $A$ by $A_1$, after which,
we may assume that $\alpha_0$ (resp. $p$) is regular
in $\bE$ (resp. in $A$). In this case, for every $c\leq 1$
we have a commutative diagram of sheaves :
$$
\xymatrix{ \cQ^c_{U_\bE} \ar[r]^-{\bar\psi{}^c_1} \ar[d] &
\phi_*\cQ^c_{U_A} \ar[d] \\
\cO_{X_\bE}/\alpha_0^c\cO_{X_\bE} \ar[r] &
\phi_*(\cO_{X_A}/p^c\cO_{X_A})
}$$
whose left (resp. right) vertical arrow is the isomorphism
induced by scalar multiplication by $\alpha_0^c$ (resp. by
$q^c$) and whose bottom horizontal arrow is the isomorphism
induced by $\phi_0^\flat$ (notation of \eqref{subsec_phi-flat-map}).
The claim follows.
\end{pfclaim}

\begin{claim}\label{cl_we-may-assume}
We may assume that the natural maps $\cO_{X_\bE}\to\bar\cO_{X_\bE}$
and $\cO_{X_A}\to\bar\cO_{X_A}$ are isomorphisms.
\end{claim}
\begin{pfclaim} Assumption (a) implies that
$Z_\bE\subset X_{\bE_2}$ and $Z_A\subset X_{A_2}$, therefore
$$
\cJ_\bE:=\Ker\,(\cO_{\!X_\bE}\to\bar\cO_{\!X_\bE})=
\underline\Gamma_{Z_\bE}\cO^\tt_{\!X_\bE}
\qquad
\cJ_A:=\Ker\,(\cO_{\!X_A}\to\bar\cO_{\!X_A})=
\underline\Gamma_{Z_A}\cO^\tt_{\!X_A}
$$
(notation of \eqref{subsec_phi-flat-map}), and lemma
\ref{lem_phi-flat-map}(ii,iii) implies that $\phi^\flat$
restricts to an isomorphism of abelian sheaves
\set\begin{equation}\label{eq_really-needed}
\cJ_\bE\isom\phi_*\cJ_A.
\end{equation}
Notice that $J:=\Gamma(X_\bE,\cJ_\bE)$ is an ideal of both
$\bE$ and $\bE_2$, and moreover $\Phi_\bE(J)=J$. Indeed, $J$
consists of all $x\in\bE$ whose images vanish in $\bE_f$, for
every $f\in\bE$ such that $\Spec\,\bE_f\subset U_\bE$; since
$\bE_f$ is still a perfect $\F_p$-algebra, it follows
easily that $x\in J$ if and only if $x^p\in J$, whence the
assertion. Thus, both $\bar\bE:=\bE/J$ and $\bar\bE_2:=\bE_2/J$
are still perfect $\F_p$-algebras, and are even perfectoid,
with the quotient topologies induced by the projections
$\bE\to\bar\bE$ and $\bE_2\to\bar\bE_2$. Moreover, it is
easily seen that the diagram $\cE$ of
\eqref{subsec_decompose-torsion} induces a natural identification
$$
\bar\bE\isom\bE_1\times_{\bE_3}\bar\bE_2.
$$
Let $h_\bE:X_{\bar\bE}:=\Spec\,\bar\bE\to X_\bE$ be the closed
immersion,
$j_{\bar\bE}:U_{\bar \bE}:=U_\bE\cap X_{\bar\bE}\to X_{\bar\bE}$
the open immersion, and consider the corresponding
quasi-coherent $\cO_{X_{\bar\bE}}$-module $\cQ^c_{U_{\bar\bE}}$ as
in \eqref{subsec_large-open-subset}.
Since $U_\bE$ is quasi-compact, $\cJ_\bE$ is a quasi-coherent
$\cO_{X_\bE}$-module (see \eqref{subsec_case-of-schs}), so
$\cJ_\bE=\Ker\,(\cO_{X_\bE}\to h_{\bE*}\cO_{X_{\bar\bE}})$; it follows
that the induced map $\cO_{X_{\bar\bE}}\to j_{\bar\bE*}\cO_{U_{\bar\bE}}$
is a monomorphism, and $h_\bE$ restricts to an isomorphism of
open subschemes $U_{\bar\bE}\isom U_\bE$, so the resulting
map $j_{\bE*}\cO_{U_\bE}\to h_{\bE*}j_{\bar\bE*}\cO_{U_{\bar\bE}}$
is an isomorphism, and we get a natural identification
$$
\omega_\bE:\cQ^c_{U_\bE}\isom h_{\bE*}\cQ^c_{U_{\bar\bE}}.
$$
Next, consider the cartesian diagram $A(\cE)$ of
\eqref{subsec_decompose-torsion}, and recall that
the bottom horizontal arrow $A_2\to A_3$ of $A(\cE)$
is naturally isomorphic to the bottom horizontal arrow
$\bE_2\to\bE_3$ of $\cE$ (remark
\ref{subsec_decompose-torsion}(iii)). There follows
a commutative diagram
$$
\xymatrix{
A \ar[r] \ar[d] & \bar A \ar[r]^-{\bar\pi_1} \ar[d]_{\bar\pi_2}
& A_1 \ar[d] \\
\bE_2 \ar[r] & \bar\bE_2 \ar[r] & \bE_3
}$$
whose two square subdiagrams are cartesian. Especially,
$J$ is naturally identified with the kernel of the
projection $A\to\bar A$; moreover, we have an
isomorphism of topological rings
$$
\bE(\bar A)\isom\bar\bE
$$
that identifies $\bE(\bar\pi_1)$ and $\bE(\bar\pi_2)$ with
the projections $\bar\bE\to\bE_1$ and $\bar\bE\to\bar\bE_2$,
and furthermore, $\bar A$ is perfectoid (proposition
\ref{prop_construct-new-perfs}). Let also
$$
h_A:X_{\bar A}:=\Spec\,\bar A\to X_A
\qquad\text{and}\qquad
j_{\bar A}:U_{\bar A}:=U_A\cap X_{\bar A}\to X_{\bar A}
$$
be respectively the closed and the open immersion, and
consider the corresponding quasi-coherent
$\cO_{X_{\bar A}}$-module $\cQ^c_{U_{\bar A}}$ defined as
in \eqref{subsec_large-open-subset}; in light of
\eqref{eq_really-needed}, we see that
$\cJ_A=\Ker\,(\cO_{\!X_A}\to h_{A*}\cO_{\!X_{\bar A}})$, and
$h_A$ restricts to an isomorphism of open subschemes
$U_{\bar A}\isom U_A$, so the resulting map
$j_{A*}\cO_{U_A}\to h_{A*}j_{\bar A*}\cO_{U_{\bar A}}$ is an
isomorphism, and we get a natural identification
$$
\omega_A:\cQ^c_{U_A}\isom h_{A*}\cQ^c_{U_{\bar A}}.
$$
Lastly, let $\bar\phi:X_{\bar A}\to X_{\bar\bE}$ be the
continuous map provided by \eqref{subsec_map-on-spec}, and
$\bar\psi{}^c_{\bar A}:\cQ^c_{U_{\bar\bE}}\to\bar\phi_*\cQ^c_{U_{\bar A}}$
the morphism of abelian sheaves attached to $\bar A$ and
$U_{\bar\bE}$ as in (i); a direct inspection yields a
commutative diagram
$$
\xymatrix{ \cQ^c_{U_\bE} \ar[rr]^-{\bar\psi{}^c} \ar[d]_{\omega_\bE}
& & \phi_*\cQ^c_{U_A} \ar[d]^{\phi_*\omega_A} \\
h_{\bE*}\cQ^c_{U_{\bar\bE}} \ar[rr]^-{h_{\bE*}\bar\psi{}^c_{\bar A}}
& & h_{\bE*}\bar\phi_*\cQ^c_{U_{\bar A}}=\phi_*h_{A*}\cQ^c_{U_{\bar A}}
}$$
whence the claim.
\end{pfclaim}

Henceforth we assume that $\cO_{X_\bE}=\bar\cO_{X_\bE}$ an
$\cO_{X_A}=\bar\cO_{X_A}$. Set $U'_\bE:=\Spec\,\bE[\alpha_0^{-1}]$,
so that $U'_A:=\phi^{-1}U'_\bE$ is the open subset
$\Spec\,A[p^{-1}]$ of $X_A$, and  let $j'_\bE:U'_\bE\to X_\bE$
and $j'_A:U'_A\to X_A$ be the open immersions. Set
also $Z'_\bE:=X_\bE\setminus U'_\bE$ and
$Z'_A:=X_A\setminus U'_A$. Notice that
$$
i_{\bE_1*}\cO_{X_{\bE_1}}=\Img(\cO_{X_\bE}\to j'_{\bE*}\cO_{U'_\bE})
\qquad\text{and}\qquad
i_{A_1*}\cO_{X_{A_1}}=\Img(\cO_{X_A}\to j'_{A*}\cO_{U'_A})
$$
and denote by
$\bar\psi{}^{\prime c}:\cQ^c_{U'_\bE}\to\cQ^c_{U'_A}$ the
map of abelian sheaves associated as in (i) with $A$ and
$U'_\bE$. Since $U'_\bE\subset U_\bE$, the natural maps
$j_{\bE*}\cO_{U_\bE}\to j'_{\bE*}\cO_{U'_\bE}$ and
$j_{A*}\cO_{U_A}\to j'_{A*}\cO_{U'_A}$ induce morphisms
$\rho_\bE:R^1\underline\Gamma_{Z_\bE}\cO_{X_\bE}\to
R^1\underline\Gamma_{Z'_\bE}\cO_{X_\bE}$ and
$\rho_A:R^1\underline\Gamma_{Z_A}\cO_{X_A}\to
R^1\underline\Gamma_{Z'_A}\cO_{X_A}$, and notice that
the image of $\rho_\bE$ (resp. of $\rho_A$) lies in
$\underline\Gamma_{Z_\bE}R^1\underline\Gamma_{Z'_\bE}\cO_{X_\bE}$
(resp. in
$\underline\Gamma_{Z_A}R^1\underline\Gamma_{Z'_A}\cO_{X_A}$);
there follows a commutative diagram
$$
\xymatrix{ \cQ^c_{U_\bE} \ar[rr]^-{\bar\psi{}^c} \ar[d]_{\rho^c_\bE}
& & \phi_*\cQ^c_{U_A} \ar[d]^{\phi_*(\rho^c_A)} \\
\underline\Gamma_{Z_\bE}\cQ^c_{U'_\bE}
\ar[rr]^-{\underline\Gamma_{Z_\bE}\bar\psi{}^{\prime c}}
& & \underline\Gamma_{Z_\bE}\phi_*\cQ^c_{U'_A}=
\phi_*\underline\Gamma_{Z_A}\cQ^c_{U'_A}
}$$
where $\rho^c_\bE$ (resp. $\rho^c_A$) is the restriction
of $\rho_\bE$ (resp. of $\rho_A$). From the snake lemma,
we obtain natural identifications
$$
\begin{aligned}
\Ker\,\rho^c_\bE\isom\,&
\Ker\,(\alpha_0^{-c}\cO_{X_\bE}\to j'_{\bE*}\cO_{U'_\bE})/
\Ker\,(\cO_{X_\bE}\to i_{\bE_1*}\cO_{X_{\bE_1}})=
(\underline\Gamma_{X_{\bE_2}}\alpha_0^{-c}\cO_{X_\bE})/\cO^\tt_{X_\bE} \\
\Ker\,\rho^c_A\isom\,&
\Ker\,(p^{-c}\cO_{X_A}\to j'_{A*}\cO_{U'_A})/
\Ker\,(\cO_{X_A}\to i_{A_1*}\cO_{X_{A_1}})=
(\underline\Gamma_{X_{A_2}}p^{-c}\cO_{X_A})/\cO^\tt_{X_A}
\end{aligned}
$$
(notation of \eqref{subsec_phi-flat-map}).

\begin{claim}\label{cl_blood-and-tears}
The map $\bar\psi{}^c$ restricts to an isomorphism of
abelian sheaves
$$
\Ker\,\rho^c_\bE\isom\Ker\,\phi_*(\rho^c_A).
$$
\end{claim}
\begin{pfclaim} Notice that
$$
\underline\Gamma_{X_{\bE_2}}\alpha^{-c}_0\cO_{\!X_\bE}=
j_{\bE*}j^*_\bE\cO^\tt_{\!X_\bE}
\qquad\text{and}\qquad
\underline\Gamma_{X_{A_2}}p^{-c}\cO_{\!X_A}=
j_{A*}j^*_A\cO^\tt_{\!X_A}.
$$
In view of lemma \ref{lem_phi-flat-map}(ii), we deduce
that $\phi^\flat$ induces an isomorphism of abelian sheaves
$$
\Ker\,\rho^c_\bE\isom
\cR':=\phi_*(\underline\Gamma_{X_{\bE_2}}\alpha^{-c}_0\cO_{\!X_\bE})/
\phi_*(\cO^\tt_{\!X_A}).
$$
Now, proposition \ref{prop_continuous-map}(ii) implies
that $\phi_*$ is an exact functor on the category of
quasi-coherent $\cO_{\!X_A}$-modules, and since both $U_A$ and
$U'_A$ are quasi-compact (see \eqref{subsec_new-place}), the
$\cO_{X_A}$-modules $j_{A*}j^*_A\cO^\tt_{\!X_A}$ and $\cO^\tt_{\!X_A}$
are quasi-coherent (\cite[Ch.I, Cor.9.2.2]{EGAI}) so $\lambda$
restricts to an isomorphism
$\cR'\isom\phi_*\Ker\,\rho^c_A=\Ker\,\phi_*(\rho^c_A)$, whence
the claim.
\end{pfclaim}

Now, claim \ref{cl_sweat-and-tears} implies that
$\underline\Gamma_{Z_\bE}\bar\psi{}^{\prime c}$ is an
isomorphism, and in light of claim \ref{cl_blood-and-tears}
we see that $\bar\psi{}^c$ is a monomorphism, so we are
reduced to showing that
$\underline\Gamma_{Z_\bE}\bar\psi{}^{\prime c}$ induces an
epimorphism $\Img\,\rho^c_\bE\to\Img\,\phi_*(\rho^c_A)$.
To this aim, let $x\in X_\bE$ be any point, and
$s_x\in\Img\,(\phi_*\rho^c_A)_x=\phi_*(\Img\,\rho^c_A)_x$
any element; by proposition \ref{prop_continuous-map}(ii),
the functor $\phi_*$ is exact on quasi-coherent
$\cO_{X_A}$-modules, so both maps
$$
\phi_*(p^{-c}\cO_{X_A})_x\to\phi_*(\cQ^c_{U_A})_x\to
\phi_*(\Img\,\rho^c_A)_x
$$
are surjective, hence we may find an open neighborhood
$\Omega\subset X_\bE$ of $x$ such that $s_x$ lifts
to a section $s\in\cO_{X_A}(U_A\cap\phi^{-1}\Omega)$
and $q^c\cdot s$ extends to a section of
$\cO_{X_A}(\phi^{-1}\Omega)$ (where $q$ is as in
\eqref{subsec_large-open-subset}). Denote by
$t_x\in(\underline\Gamma_{Z_\bE}\cQ^c_{U'_\bE})_x$
the preimage of $s_x$. Notice that
$$
\underline\Gamma_{Z_\bE}R^1\underline\Gamma_{Z'_\bE}\cO_{X_\bE}=
(j_{\bE*}j^*_{\bE}(i_{\bE_1*}\cO_{X_{\bE_1}}))/i_{\bE_1*}\cO_{X_{\bE_1}}
$$
so we may replace $\Omega$ by a smaller open neighborhood
of $x$, and assume that $t_x$ is the image of a section
$\bar t\in i_{\bE_1*}\cO_{X_{\bE_1}}(U_\bE\cap\Omega)$. Denote
by $\bar s\in i_{A_1*}\cO_{X_{A_1}}(U_A\cap\phi^{-1}\Omega)$
the image of $s$, and set
$$
\bar u:=i_{\bE_1*}\phi^\flat_1(\bar t)-\bar s\in
j_{A_*}j^*_Ai_{A_1*}\cO_{X_{A_1}}(\phi^{-1}\Omega).
$$
By construction, the image of $\bar u_x$ vanishes in
$(\phi_*R^1\underline\Gamma_{Z'_A}\cO_{X_A})_x$, thus
$\bar u_x\in\phi_*(i_{A_1*}\cO_{X_{A_1}})_x$. By proposition
\ref{prop_continuous-map}(ii), the natural map
$\phi_*(\cO_{X_A})_x\to\phi_*(i_{A_1*}\cO_{X_{A_1}})_x$
is surjective, so we may replace $\Omega$ by a smaller
open neighborhood of $x$, and assume that $\bar u$ is the
image of a local section $u\in\cO_{X_A}(\phi^{-1}\Omega)$.
However, $s$ and $s':=s+u_{|U_A\cap\phi^{-1}\Omega}$ have the
same image in $\phi_*(\cQ^c_{U_A})_x$, so we may replace $s$
by $s'$, and assume that $i_{\bE_1*}\phi^\flat_1(\bar t)=\bar s$.
In this situation, remark \ref{rem_phi-is-natural}(ii)
implies that $\bar t$ lifts to a section
$t\in\cO_{X_\bE}(U_\bE\cap\Omega)$ such that $\phi^\flat(t)=s$,
and we denote by
$v\in R^1\underline\Gamma_{Z_\bE}\cO_{X_\bE}(\Omega)$ the
image of $t$. Clearly, the image of $v$ in
$(\underline\Gamma_{Z_\bE}R^1\underline\Gamma_{Z'_\bE}\cO_{X_\bE})_x$
agrees with $t_x$, hence it remains only to check that
$v_x\in\cQ_{U_\bE,x}^c$. But notice that $\alpha_0^c\cdot t_x=0$,
so $\alpha_0^c\cdot v_x\in(j_{\bE*}j^*_\bE\cO^\tt_{X_\bE})_x$.
Pick any $\eps\in\N[1/p]$ such that
$0<\eps\leq\min(p/(p-1)-c,1)$; by virtue of corollary
\ref{cor_perf-are-reduced}(i), it follows that
$\alpha_0^{c+\eps}\cdot v_x=0$, {\em i.e.}
$\alpha_0^c\cdot v_x\in\cQ^\eps_{U_\bE,x}$. Lastly, from (i)
we get
$$
\bar\psi{}^\eps(\alpha_0^cv_x)=
q^c\cdot\bar\psi{}^{c+\eps}(v_x)
$$
and $\bar\psi{}^{c+\eps}(v_x)$ is the image of $s$ in
$\phi_*(\cQ^{c+\eps}_{U_A})_x$. By construction, the
image of $q^c\cdot s$ vanishes in $\phi_*(\cQ^{c+\eps}_{U_A})_x$,
and therefore $\alpha_0^c\cdot v_x=0$, since we have already
remarked that $\bar\psi{}^\eps$ is a monomorphism. The
assertion follows.
\end{proof}

\subsection{Perfectoid quasi-affinoid rings}
We wish now to merge the theory of perfectoid rings with
that of affinoid rings of section \ref{sec_affinoid-rings}.
First, let us make the following :

\begin{definition}\label{def_perfectoid-qaff}
(i)\ \
We say that a quasi-affinoid ring $(A,A^+,U)$ is
{\em perfectoid} if $A$ is a perfectoid ring.

(ii)\ \
We say that a quasi-affinoid scheme $\underline X$ is
{\em perfectoid}, if there exists a perfectoid
quasi-affinoid ring $\underline A$ and an isomorphism
$\underline X\isom\sSpec\,\underline A$.

(iii)\ \
We denote by
$$
\mathsf{q.Afd.Ring}_\mathrm{perf}
\qquad\text{and}\qquad
\mathsf{q.Afd.Sch}_\mathrm{perf}
$$
the subcategories of $\mathsf{q.Afd.Ring}$ and
respectively $\mathsf{q.Afd.Sch}$ whose objects
are the perfectoid quasi-affinoid rings (resp.
schemes), and whose morphisms are the f-adic
morphisms of quasi-affinoid rings (resp. schemes).
\end{definition}

\sset\subsubsection{}\label{subsec_general-case}
Clearly we have a well defined functor
$$
\mathsf{q.Afd.Ring}_\mathrm{perf}\to
\mathsf{q.Afd.Sch}^o_\mathrm{perf}
\qquad
\underline A\mapsto\sSpec\,\underline A
$$
and we will construct a right adjoint. To this aim, we
need a few preliminaries : we consider any perfectoid
quasi-affinoid ring $\underline A:=(A,A^+,U_A)$ and set
$X_A:=\Spec\,A$, $Z_A:=X_A\setminus U_A$. Let also
$\bE:=\bE(A)$, define the continuous map
$\phi:X_A\to X_\bE:=\Spec\,\bE$ as in
\eqref{subsec_map-on-spec}, and set $Z_\bE:=\phi(Z_A)$
and $U_\bE:=X_\bE\setminus Z_\bE$. Then it is easily seen
that $U_\bE$ contains the analytic locus of $X_\bE$, so
any choice of ring of integral elements $\bE^+\subset\bE$
will give another perfectoid quasi-affinoid ring
$\underline\bE:=(\bE,\bE^+,U_\bE)$, as well as the
attached perfectoid quasi-affinoid schemes
$$
(U_A,\cT_{U_A},A^+_U):=\sSpec\,\underline A
\qquad
(U_\bE,\cT_{U_\bE},\bE_U^+):=\sSpec\,\underline\bE.
$$
More precisely, let $\underline\alpha:=(\alpha_n~|~n\in\N)$
be any distinguished element in $\Ker\,u_A$; since $\phi$
maps $\Spec\,A/pA$ bijectively onto $\Spec\,\bE/\alpha_0\bE$,
and maps $\Spec\,A[p^{-1}]$ onto $\Spec\,\bE[\alpha^{-1}_0]$,
it is easily seen that $U_\bE$ is the unique open subset of
$X_\bE$ containing the analytic locus, and such that
\set\begin{equation}\label{eq_better}
\phi^{-1}(U_\bE)=U_A.
\end{equation}
Recall that, according to lemma \ref{lem_topology-on-opens},
the topologies $\cT_{U_A}$ on $A_U:=\cO_{\!X_A}(U_A)$ and
$\cT_{U_\bE}$ on $\bE_U:=\cO_{\!X_\bE}(U_\bE)$ are f-adic,
complete and separated, and the restriction maps $A\to A_U$
and $\bE\to\bE_U$ are open. Especially the image $\bar A$
of $A$ (resp. $\bar\bE$ of $\bE$) is open in $A_U$ (resp.
in $\bE_U$), and we endow this subring with the quotient
topology induced by the projection $\pi_A:A\to\bar A$
(resp. $\pi_\bE:\bE\to\bar\bE$). Moreover, we have a
continuous map of (multiplicative) monoids
$$
\phi^\flat_U:\bE_U\to A_U
$$
(proposition \ref{prop_new-formula}(i)). Let also
$$
j^\circ_\bE:\bar\bE\to\bE_U^\circ
\qquad\text{and}\qquad
j^\circ_A:\bar A\to A_U^\circ
$$
be the inclusion maps, and we endow $A^\circ_U$ (resp.
$\bE^\circ_U$) with the topology induced by the inclusion
into $A_U$ (resp. into $\bE_U$), so that $j^\circ_A$ and
$j^\circ_\bE$ are open maps.

\begin{proposition}\label{prop_general-case}
In the situation of \eqref{subsec_general-case}, we have :
\begin{enumerate}
\item
For every $x\in A_U$, the following conditions are equivalent :
\begin{enumerate}
\item
$x$ is power bounded in $A_U$.
\item
The subset $\{x^{p^n}~|~n\in\N\}$ is bounded in $A_U$.
\item
The image of $x$ in $R^1\Gamma_{\!Z_A}\cO_{X_A}$ is annihilated
by $A^{\circ\circ}$.
\item
$x\cdot A^{\circ\circ}_U\subset A^{\circ\circ}_U$.
\item
$v(x)\leq 1$ for every analytic rank one valuation
$v\in\Spa\,\underline A$.
\end{enumerate}
\item
$\bar A$ and\/ $\bar\bE$ are perfectoid rings, and there exists
a unique isomorphism
$$
\bar\omega:\bar\bE\isom\bE(\bar A)
\qquad\text{such that}\qquad
\bar\omega\circ\pi_\bE=\bE(\pi_A).
$$
\item
The maps $\pi_A$ and $\pi_\bE$ restrict to bijections
$\bE^{\circ\circ}\isom\bE_U^{\circ\circ}$ and
$A^{\circ\circ}\isom A_U^{\circ\circ}$.
\item
The map $\phi^\flat_U$ restricts to a continuous map of
topological monoids
$$
\phi_U^{\flat\circ}:\bE_U^\circ\to A_U^\circ
$$
which induces an isomorphism of perfect $\F_p$-algebras
$$
\bar\phi{}^{\flat\circ}_U:
\bE_U^\circ/\bE_U^{\circ\circ}\isom A^\circ_U/A_U^{\circ\circ}.
$$
\end{enumerate}
\end{proposition}
\begin{proof} We consider first the following special case :

\begin{claim}\label{cl_power-bounded}
In the situation of \eqref{subsec_general-case}, let
$\beta\in\bE$ be any regular element such that
$\alpha_0\in\beta\bE$, set $b:=\bar u_A(\beta)\in A$,
and suppose moreover that $bA$ is an ideal of adic
definition for $A$, and $U_A$ is the analytic locus
of $X_A$. Then conditions (i.a)--(i.d) are equivalent
in this case, and
$$
A^{\circ\circ}=A_U^{\circ\circ}=
\bigcup_{\lambda\in\N[1/p]\setminus\{0\}}b^\lambda A
\qquad
\text{with $b^\lambda:=\bar u_A(\beta^\lambda)$ for every
$\lambda\in\N[1/p]$}.
$$
\end{claim}
\begin{pfclaim} First, notice that under the assumptions
of the claim, $p\in bA$ (lemma \ref{lem_was-third-cond}(iii))
and $b^\lambda$ is a regular element of $A$ for every
$\lambda\in\N[1/p]$ (proposition \ref{prop_back-and-forth-sec}).
Especially, the localization map $A\to A[b^{-1}]$ is injective,
and $U_A=\Spec\,A[b^{-1}]$, so that $A_U=A[b^{-1}]$,
$R^1\Gamma_{\!Z_A}\cO_{\!X_A}=A[b^{-1}]/A$, and condition (i.c) implies :
\begin{itemize}
\item[(i.f)]
$b^\lambda x\in A$ for every $\lambda\in\N[1/p]\setminus\{0\}$.
\end{itemize}
We show first the equivalence of (i.a), (i.b) and (i.f).
Indeed, obviously (i.f)$\Rightarrow$(i.a)$\Rightarrow$(i.b),
so we may assume that (i.b) holds, and we set
$$
\rho:=\inf\{\gamma\in\N[1/p]~|~
b^\gamma x^{p^n}\in A\ \text{for every $n\in\N$}\}.
$$
Notice that (i.b) implies that $\rho<+\infty$. We claim
that $\rho=0$. Indeed suppose, by way of contradiction,
that $\rho>0$ and pick any $\gamma\in\N[1/p]$ such that
$\rho>\gamma>\max(\rho/p,\rho-1/p)$.
Then $b^{1/p}b^\gamma x^{p^n}\in A$, and hence
$b^\gamma x^{p^n}\in b^{-1/p}A$ for every $n\in\N$; on the
other hand, we also see that $b^{p\gamma}x^{p^n}\in A$ for
every $n\in\N$. Now, by corollary \ref{cor_jackob}, the
Frobenius endomorphism of $A/pA$ induces an isomorphism
$A/b^{1/p}A\to A/bA$; since $b$ is regular, there follows
a bijective map
$$
b^{-1/p}A/A\isom b^{-1}A/A
\qquad
(x\mod A)\mapsto(x^p\mod A).
$$
We deduce that $b^\gamma x^{p^n}\in A$ for every $n\in\N$,
{\em i.e.} $\gamma\leq\rho$, which is absurd. This shows
that (i.f) holds, whence the assertion.
Next, say that $x\in A[b^{-1}]^{\circ\circ}$; since $bA$
is open in $A[b^{-1}]$, we must have $x^{p^n}\in bA$ for
some $n\in\N$, so if set $y:=b^{-1/p^n}x$, we get $y^{p^n}\in A$,
and therefore the subset $\{y^{p^n}~|~n\in\N\}$ is bounded
in $A_U$. Pick any integer $r>n$; by the foregoing, it
follows that $b^{1/p^r}y\in A$, and finally
$x\in b^{1/p^n-1/p^r}A$. This also implies the equivalence
of (i.c),(i.d) and (i.f), so the proof of the claim is
complete.
\end{pfclaim}

For the general case, pick a finite system
$\beta_\bullet:=(\beta_0,\dots,\beta_k)$ of elements of
$\bE$ with $\beta_0=\alpha_0$, and such that $\beta_\bullet$
generates an ideal $\cI$ of adic definition for $\bE$. Then
the ideal $I\subset A$ generated by
$(\bar u_A(\beta_0),\dots,\bar u_A(\beta_k))$
is of adic definition for $A$ (see remark
\ref{rem_nice-topology}(iii)). We notice :

\begin{claim}\label{cl_identify-top-nilp}
The ideal $A^{\circ\circ}$ is generated by
$(\bar u_A(\beta_i^\lambda)~|~i=0,\dots,k,\ 
\lambda\in\N[1/p]\setminus\{0\})$.
\end{claim}
\begin{pfclaim} Clearly
$\bar u_A(\beta_i^\lambda)\in A^{\circ\circ}$ for every
$i=0,\dots,k$ and every strictly positive $\lambda\in\N[1/p]$.
Conversely, let $c\in A^{\circ\circ}$ be any element; then
$c^{p^n}\in I$ for every sufficiently large $n\in\N$.
By remark \ref{rem_nice-topology}(iii) we may find
$\gamma\in\bE$ such that $\bar u_A(\gamma)-c\in pA$,
and the image of $\gamma^{p^n}$ in $\bE/\alpha_0\bE$ lies
in $\cI/\alpha_0\bE$, so $\gamma^{p^n}\in\cI$, hence
$\gamma\in\cI^{\La 1/p^n\Ra}\bE$, and therefore
$c\in\cI^{\La 1/p^n\Ra}A+pA=\cI^{\La 1/p^n\Ra}A$ (notation
of remark \ref{rem_Witt-are-f-adic}(i)), whence the
claim.
\end{pfclaim}

(i.c)$\Rightarrow$(i.a): By assumption,
$\bar u_A(\beta^\lambda_i)\cdot x$ lies in $\bar A$,
for every $i=0,\dots,k$ and every strictly positive
$\lambda\in\N[1/p]$; let $n>0$ be any integer, and pick
$t\in\N$ such that $p^t\geq n$; then
$\bar u_A(\beta^{n/p^t}_i)\cdot x^n=
(\bar u_A(\beta^{1/p^t})\cdot x)^n$ lies as well in
$\bar A$, and we conclude that
$\bar u_A(\beta_i)\cdot x^n\in\bar A$ for every
$n\in\N$, whence the assertion.

Obviously, (i.a)$\Rightarrow$(i.b). To show that
(i.b)$\Rightarrow$(i.c), let $\cT_{A,p}$ be the $p$-adic
topology on $A$, and $\cT_{U,p}$ the unique f-adic topology
on $A_U$ such that the restriction map
$(A,\cT_{A,p})\to(A_U,\cT_{U,p})$ is open (lemma
\ref{lem_topology-on-opens}(i)). Suppose first that the
topology of $A$ agrees with $\cT_{A,p}$, and let
$$
A_1:=A/\Ann_A(p)
\qquad
Z':=\Spec\,A/pA
\qquad
U':=X_A\setminus Z'.
$$
Denote as well by $x_{|U'}\in\cO_{X_A}(U')=A[1/p]$ the image
of $x$. We recall that $p$ is regular in $A_1$ (remark
\ref{rem_decompose-torsion}(ii)), and we endow $A_1$
with its $p$-adic topology $\cT_{A_1,p}$, and $A[1/p]$ with
the unique f-adic topology $\cT_{U',p}$ such that the
inclusion map $(A_1,\cT_{A_1,p})\to(A[1/p],\cT_{U',p})$ is open,
so that the subset $\{x^{p^n}_{|U'}~|~n\in\N\}$ is bounded
in $A[1/p]$. Set $q:=\bar u_A(\alpha_0)$; by claim
\ref{cl_power-bounded} and lemma
\ref{lem_was-third-cond}(iii), it then follows that
$q^\lambda\cdot x_{|U'}$ lies in $A_1$ for every
$\lambda\in\N[1/p]\setminus\{0\}$.
Denote also by $\bar x\in R^1\Gamma_Z\cO_{X_A}$ the
image of $x$; we deduce that the image of
$q^\lambda\cdot\bar x$ vanishes in
$R^1\Gamma_{Z'}\cO_{X_A}$, for every such $\lambda$.
To conclude in this case, it suffices to remark :

\begin{claim} The kernel of the natural map
$R^1\Gamma_{Z_A}\cO_{X_A}\to R^1\Gamma_{Z'}\cO_{X_A}$
is annihilated by $q^\lambda$ for every
$\lambda\in\N[1/p]\setminus\{0\}$.
\end{claim}
\begin{pfclaim} By the snake lemma, this kernel is a
quotient of $\Gamma_{\!Z'}j_{A*}\cO_{U_A}$, where $j_A:U_A\to X_A$
is the open immersion. However,
$\Gamma_{\!Z'}j_*\cO_{U_A}=H^0(U_A,\underline\Gamma_{Z'}\cO_{X_A})$.
Since $\underline\Gamma_{Z'}\cO_{X_A}$ is a quasi-coherent
$\cO_{X_A}$-module (lemma \ref{lem_without-cohereur}(i)),
we are reduced to checking that $\Ann_A(p)$ is annihilated
by $q^\lambda$ for every $\lambda\in\N[1/p]\setminus\{0\}$.
The latter holds by corollary \ref{cor_perf-are-reduced}(i).
\end{pfclaim}

Next, let $A$ be any perfectoid ring; let $j_\bE:U_\bE\to X_\bE$
be the open immersion, and for every $c\in\N[1/p]$ consider
the $\cO_{X_\bE}$-module $\cQ^c_{U_\bE}$ and the $\cO_{X_A}$-module
$\cQ^c_{U_A}$ as in \eqref{subsec_large-open-subset}. We set
$$
\cQ^0_{U_\bE}:=\bigcap_{c>0}\cQ^c_{U_\bE}
\qquad\text{and}\qquad
\cQ^0_{U_A}:=\bigcap_{c>0}\cQ^c_{U_A}.
$$
Condition (i.b) implies that the set $\{x^{p^n}~|~n\in\N\}$ is
bounded also for the topology $\cT_{U,p}$ of $A_U$, and since
$(A,\cT_{A,p})$ is perfectoid (proposition
\ref{prop_change-topol}(i)), the foregoing case shows that
$\bar x$ lies in $\cQ^0_{U_A}(X_A)$.
By proposition \ref{prop_Q-c-tricky}(iv), the map
$\phi^\flat$ of lemma \ref{lem_phi-flat-map}(i) induces
an isomorphism
$$
\bar\psi{}^0:\cQ^0_{U_\bE}\isom\phi_*\cQ^0_{U_A}
$$
of abelian sheaves, fitting into a commutative diagram
\set\begin{equation}\label{eq_psi-and-bp}
{\diagram
\cQ^0_{U_\bE} \ar[r]^-{\bar\psi{}^0} \ar[d]_{\bar\bp_\bE} &
\phi_*\cQ^0_{U_A} \ar[d]^{\phi_*(\bar\bp_A)} \\
\cQ^0_{U_\bE} \ar[r]^-{\bar\psi{}^0} & \phi_*\cQ^0_{U_A}
\enddiagram}
\end{equation}
where $\bar\bp_\bE$ (resp. $\bar\bp_A$) is induced by
the $p$-Frobenius endomorphism of the sheaf of monoids
$j_{\bE*}\cO_{U_\bE}$ (resp. $j_{A*}\cO_{U_A}$). With this
notation, condition (i.b) says that for some $n\in\N$
we have
\set\begin{equation}\label{eq_summer}
\bar u_A(\beta_i^n)\cdot\bar\bp_A^r(\bar x)=0
\qquad
\text{for every $i=0,\dots,k$ and every $r\in\N$}.
\end{equation}
Let $\bar\gamma\in\cQ^0_{U_\bE}(X_\bE)$ be the preimage of
$\bar x$; then $\beta_i^\lambda\cdot\bar\gamma$ is the
preimage of the class of
$\bar u_A(\beta_i^\lambda)\cdot x$ in $\cQ^0_{U_A}(X_A)$,
for every $i=0,\dots,k$ and every $\lambda\in\N[1/p]$.
By virtue of claim \ref{cl_identify-top-nilp}, we are then
reduced to checking that $\beta_i^\lambda\cdot\bar\gamma=0$
for every such $i$ and $\lambda$, and \eqref{eq_summer}
yields already
$$
\beta_i^n\cdot\bar\bp^r_\bE(\bar\gamma)=0
\qquad
\text{for every $i=0,\dots,k$ and every $r\in\N$}.
$$
Now, we may write $\lambda=t\cdot p^{-r}$ for some integers
$t\geq n$ and $r\in\N$; we deduce
$$
\bar\bp^r_\bE(\beta_i^\lambda\cdot\bar\gamma)=
\beta_i^t\cdot\bar\bp^r_\bE(\bar\gamma_\bE)=0.
$$
Lastly, since $\bE$ is perfect, it is easily seen that
$\bar\bp_\bE$ is an isomorphism, whence the assertion.

(ii): This was already remarked in the proof of claim
\ref{cl_we-may-assume}. Let us add the following :

\begin{claim}\label{cl_nilp-in-bar-A}
The map $\pi_\bE$ and $\pi_A$ restrict to bijections
$\bE^{\circ\circ}\isom\bar\bE{}^{\circ\circ}$ and
$A^{\circ\circ}\isom\bar A{}^{\circ\circ}$.
\end{claim}
\begin{pfclaim} First, notice that the image of $\beta_\bullet$
in $\bar\bE$ generates an ideal of adic definition for $\bar\bE$;
in view of (ii) and claim \ref{cl_identify-top-nilp}, we
deduce immediately that $\bar A{}^{\circ\circ}=\pi_A(A^{\circ\circ})$
and $\bar\bE{}^{\circ\circ}=\pi_\bE(\bE^{\circ\circ})$. Next, if
$I$ is any ideal of adic definition for $A$, clearly we
have $I\cdot\Ker\,\pi_A=0$, whence $I\cap\Ker\,\pi_A=0$, since
$A$ is reduced (corollary \ref{cor_perf-are-reduced}); but
$A^{\circ\circ}$ is the union of all ideals of adic definition
for $A$, so we conclude that $A^{\circ\circ}\cap\Ker\,\pi_A=0$.
Likewise, we see that $\bE^{\circ\circ}\cap\Ker\,\pi_\bE=0$,
whence the claim.
\end{pfclaim}

(iii): Let $x\in A_U^{\circ\circ}$ be any element; denote by
$\bar x\in R^1\Gamma_{\!Z_A}\cO_{X_A}$ the image of $x$; from
(i) we see that $\bar x\in Q_A:=\cQ^0_{U_A}(X_A)$, and since
$\bar A$ is open in $A_U$, there exists $n\in\N$ such that
$x^{p^n}\in\bar A$, therefore $\bar\bp_A^n(\bar x)=0$. On
the other hand, we have already remarked that the maps
$\bar\psi{}^0$ and $\bar\bp_\bE$ appearing in
\eqref{eq_psi-and-bp} are isomorphisms, so the same holds
for $\bar\bp_A$, and consequently $\bar x=0$, {\em i.e.}
$x\in\bar A$; then the assertion follows from claim
\ref{cl_nilp-in-bar-A}. Next, from claim
\ref{cl_identify-top-nilp} it is easily seen that the
isomorphism $\bE/\alpha_0\bE\isom A/pA$ of remark
\ref{rem_nice-topology}(ii) restricts to a natural
identification
\set\begin{equation}\label{eq_identify-circ-circ}
\bE^{\circ\circ}/\alpha_0\bE\isom A^{\circ\circ}/pA
\end{equation}
therefore
$Z_\bE\subset X_\bE\setminus\Spec\,\bE/\bE^{\circ\circ}$,
so the foregoing applies as well to the perfectoid
ring $\bE$ and its closed subset $Z_\bE$, and the proof
of (iii) is thus complete.

Now, taking into account (iii), we easily see that
(i.d)$\Rightarrow$(i.c). Conversely, if (i.c) holds,
we get $x\cdot A^{\circ\circ}\cdot A^{\circ\circ}\subset A^{\circ\circ}_U$,
but claim \ref{cl_identify-top-nilp} implies that
$A^{\circ\circ}\cdot A^{\circ\circ}=A^{\circ\circ}$, whence (i.d).
Lastly, from remark \ref{rem_about-Spa}(iv) it is clear
that (i.a)$\Rightarrow$(i.e). Conversely, suppose (i.e)
holds; then we have as well $v(xa)<1$ for every
$a\in A^{\circ\circ}_U$ and every rank one analytic valuation
$v\in\Spa\,\underline A$. Taking into account lemma
\ref{lem_Cont-A}(v), it follows that $v(xa)<1$ for every
$a\in A^{\circ\circ}_U$ and every $v\in(\Spa\,\underline A)_\mathrm{a}$.
On the other hand, we have $v(xa)=0$ for every
$v\in(\Spa\,\underline A)_\mathrm{na}$, and we conclude that
$xa\in A^{\circ\circ}_U$, by corollary \ref{cor_corcor}(ii).
This shows that (i.e)$\Rightarrow$(i.d), and completes
the proof of (i).

(iv): Set also $Q_\bE:=\cQ^0_{U_\bE}(X_\bE)$; in light of
\eqref{eq_identify-circ-circ}, it is easily seen that
the isomorphism $\bar\psi{}^0$ restricts to a natural
identification
$$
\Ann_{Q_\bE}(\bE^{\circ\circ})\isom\Ann_{Q_A}(A^{\circ\circ})
$$
and taking into account (i) we deduce that
$\phi^\flat_U(\bE^\circ_U)\subset A_U^\circ$, so the sought
map $\phi_U^{\flat\circ}$ is well defined, and it is continuous,
by proposition \ref{prop_new-formula}(i); combining with
(iii), we also see that $\phi^{\flat\circ}_U$ induces
bijections
$$
\bE_U^\circ/\bar\bE\isom A_U^\circ/\bar A
\qquad
\bar\bE/\bE_U^{\circ\circ}\isom\bar A/A_U^{\circ\circ}.
$$
On the other hand, we remark :

\begin{claim}\label{cl_slightly-odd-placing}
(i)\ \
For every $x,y\in\bE^\circ_U$ we have
$\phi^\flat_U(x+y)-\phi^\flat_U(x)-\phi^\flat_U(y)\in A_U^{\circ\circ}$.
\begin{enumerate}
\addenu
\item
$\phi^{\flat\circ}_U(\bE_U^{\circ\circ})\subset A_U^{\circ\circ}$.
\end{enumerate}
\end{claim}
\begin{pfclaim}(i): By proposition \ref{prop_new-formula}(ii),
the difference in the claim can be written as a $p$-adically
convergent series $\sum_{n\in\N\setminus\{0\}}p^n\cdot c_n$, where
each summand $c_n$ is a finite $\Z_p$-linear combination
of terms of the form $\phi^\flat_U(x^\lambda y^{1-\lambda})$,
where $\lambda,1-\lambda\in\N[1/p]$. It is easily seen
that $x^\lambda y^{1-\lambda}\in\bE^\circ_U$ for every such
$\lambda$, and obviously $p\in A^{\circ\circ}_U$, whence
the contention.

(ii) is obvious, since $\phi^{\flat\circ}_U$ is a continuous
morphism of multiplicative monoids.
\end{pfclaim}

From claim \ref{cl_slightly-odd-placing} it follows
that $\phi^{\flat\circ}_U$ descends to a map
$\bar\phi{}^{\flat\circ}_U:
\bE_U^\circ/\bE_U^{\circ\circ}\to A^\circ_U/A_U^{\circ\circ}$
of abelian groups, which must then be bijective, by the
foregoing. By remark \ref{rem_someth-on-bdd-in-Z-lin}(iv),
both the source and target of $\bar\phi{}^{\flat\circ}_U$ are
endowed with natural quotient ring structures, and it is
easily seen that $\bE_U^\circ/\bE_U^{\circ\circ}$ is a perfect
$\F_p$-algebra (details left to the reader); to conclude
the proof it then suffices to notice that, by construction,
$\bar\phi{}^{\flat\circ}_U$ is also a morphism of multiplicative
monoids.
\end{proof}

\begin{theorem}\label{th_int-subrings-perfectoid}
In the situation of \eqref{subsec_general-case}, the
following holds :
\begin{enumerate}
\item
$\bE^\circ_U$ and $A^\circ_U$ are both perfectoid, and
there exists an isomorphism of topological rings
$$
\omega^\circ:\bE^\circ_U\isom\bE(A^\circ_U)
\qquad\text{such that}\qquad
\omega^\circ\circ j^\circ_\bE=\bE(j^\circ_A).
$$
\item
The rule $D\mapsto\bE(D)$ establishes a bijection from
the set $\cS_A$ of open integrally closed subrings of
$A_U^\circ$ to the set $\cS_\bE$ of open integrally
closed subrings of\/ $\bE_U^\circ$.
\item
Moreover, every open integrally closed subring of
$A_U^\circ$ is perfectoid (for the topology induced
by $A_U^\circ$).
\end{enumerate}
\end{theorem}
\begin{proof}(i): First, notice that $\bE_U^{\circ\circ}$
(resp. $\bA^{\circ\circ}_U$) is a bounded open ideal of
$\bE^\circ_U$ (resp. of $A^\circ_U$), by virtue of proposition
\ref{prop_general-case}(iii); it follows that
$\bE^\circ_U$ and $A^\circ_U$ are bounded, so their
topologies are both adic and f-adic (corollary
\ref{cor_f-adics}(iii)). Moreover, since $\bE_U$ and
$A_U$ are complete and separated (lemma
\ref{lem_topology-on-opens}(ii)), the same holds for
$\bE^\circ_U$ and $A^\circ_U$. Furthermore, since $\bE_U$
is a perfect $\F_p$-algebra, proposition
\ref{prop_general-case}(i) easily implies that the same
holds for $\bE^\circ_U$, so the latter is perfectoid (example
\ref{ex_perfectoid}(i)).
Next, let $I\subset A$ be any ideal of definition;
proposition \ref{prop_general-case}(iii) implies that
$I_U:=IA_U^\circ\subset A^{\circ\circ}$, so $I_U$ is a finitely
generated ideal of adic definition for $A^\circ_U$ (corollary
\ref{cor_f-adics}(iii)), and clearly $p\in I_U^2$.
Furthermore, $I_U\subset\bar A$ (by claim \ref{cl_nilp-in-bar-A}),
and the Frobenius endomorphism $\Phi_{A_U/I_U}$ of $A_U/I_U$
induces surjective endomorphisms of both $\bar A/I_U$ and
$A_U/\bar A$ (proposition \ref{prop_general-case}(iv)); thus,
$\Phi_{A_U/I_U}$ is surjective, which shows that $A_U$ is a P-ring.
Next, endow $A^\circ_U/pA^\circ_U$ with the quotient topology
induced by the projection $\pi_U:A^\circ_U\to A^\circ_U/pA^\circ_U$,
and let $\phi^{\flat\circ}_U:\bE^\circ_U\to A^\circ_U$ and
$\bar\omega:\bar\bE\isom\bE(\bar A)$ be respectively
the continuous morphism of topological monoids and the
isomorphism of topological rings provided by proposition
\ref{prop_general-case}(iii,iv); in light of proposition
\ref{prop_new-formula}(ii) it is easily seen that
$\pi_U\circ\phi^{\flat\circ}_U:\bE^\circ_U\to A^\circ_U/pA^\circ_U$
is a continuous ring homomorphism fitting into a commutative
diagram
$$
\xymatrix{
\bar\bE \ar[rr]^-{\bar u_{\bar A/p\bar A}\circ\bar\omega}
\ar[d]_{j^\circ_\bE} & &
\bar A/p\bar A \ar[d]^{j^\circ_A\otimes_\Z\F_p} \\
\bE^\circ_U \ar[rr]^-{\pi_U\circ\phi^{\flat\circ}_U} & &
A^\circ_U/pA^\circ_U.
}$$
By proposition \ref{prop_lift-Witt}(iii), there follows
a continuous ring homomorphism
$v_U:W(\bE^\circ_U)\to A^\circ_U$ fitting into a commutative
diagram
$$
\xymatrix{
W(\bar\bE) \ar[rr]^-{u_{\bar A}\circ W(\bar\omega)}
\ar[d]_{W(j^\circ_\bE)} & & \bar A \ar[d]^{j^\circ_A} \\
W(\bE^\circ_U) \ar[rr]^-{v_U} & & A^\circ_U.
}$$
Since both $u_{\bar A}$ and $j^\circ_A$ are open, the same
holds for $v_U$, and since $\bar A$ is perfectoid
(proposition \ref{prop_general-case}(iii)), the kernel of
$u_{\bar A}\circ W(\bar\omega)$ is a distinguished ideal,
so it remains only to check

\begin{claim} $v_U$ induces a ring isomorphism
$$
\bar v_U:W(\bE^\circ_U)\otimes_{W(\bar\bE)}\bar A\isom A^\circ_U.
$$
\end{claim}
\begin{pfclaim} Let $\Phi_{\bE^\circ_U}$ be the Frobenius
endomorphism of $\bE^\circ_U$, and notice that
$$
\Phi_{\bE^\circ_U}(\bE^{\circ\circ}_U)=\bE^{\circ\circ}_U.
$$
Moreover, $W(\bE^{\circ\circ}_U)$ is a closed ideal of both
$W(\bar\bE)$ and $W(\bE^\circ_U)$ (remark
\ref{rem_Witt-limit}(iv)); there follows a commutative
ladder with exact rows :
$$
\cL
\quad : \quad
{\diagram
0 \ar[r] & W(\bE^{\circ\circ}_U) \ar[r] \ddouble &
W(\bar\bE) \ar[r] \ar[d] &
W(\bar\bE/\bE^{\circ\circ}_U) \ar[r] \ar[d] & 0 \\
0 \ar[r] & W(\bE^{\circ\circ}_U) \ar[r] &
W(\bE^\circ_U) \ar[r] &
W(\bE^\circ_U/\bE^{\circ\circ}_U) \ar[r] & 0.
\enddiagram}$$
Let $\underline\alpha$ be any distinguished element of
$\Ker\,u_A$; since the image of $\underline\alpha$ is a
regular element in both $W(\bE^\circ_U/\bE^{\circ\circ}_U)$
and $W(\bar\bE/\bE^{\circ\circ}_U)$, we have
$$
\Tor_1^{W(\bE)}(W(\bE^\circ_U/\bE^{\circ\circ}_U),A)=0
\qquad
\Tor_1^{W(\bE)}(W(\bar\bE/\bE^{\circ\circ}_U),A)=0
$$
so the rows of the ladder $\cL\otimes_{W(\bar\bE)}\bar A$
are still exact. Furthermore, since the topologies of
$\bar\bE/\bE^{\circ\circ}_U$ and $\bE^\circ_U/\bE^{\circ\circ}_U$
are discrete, the isomorphism $\bar\phi{}^{\flat\circ}_U$ of
proposition \ref{prop_general-case}(iv) induces identifications
$$
W(\bar\bE/\bE^{\circ\circ}_U)\otimes_{W(\bar\bE)}\bar A\isom
\bar A/A^{\circ\circ}_U
\qquad
W(\bE^\circ_U/\bE^{\circ\circ}_U)\otimes_{W(\bar\bE)}\bar A\isom
A_U^\circ/A^{\circ\circ}_U.
$$
Summing up, we obtain the commutative diagram with exact rows :
$$
\xymatrix{ 0 \ar[r] & A^{\circ\circ}_U \ar[r] \ddouble &
\bar A \ar[r] \ar[d]_w & \bar A/A^{\circ\circ}_U \ar[r] \ar[d] & 0 \\
0 \ar[r] & A^{\circ\circ}_U \ar[r] \ddouble &
W(\bE^\circ_U)\otimes_{W(\bar\bE)}\bar A \ar[r] \ar[d]_{\bar v_U} &
A^\circ_U/A^{\circ\circ}_U \ar[r] \ddouble & 0 \\
0 \ar[r] & A^{\circ\circ}_U \ar[r] &
A^\circ_U \ar[r] & A^\circ_U/A^{\circ\circ}_U \ar[r] & 0
}$$
and to conclude, it suffices to check that
$\bar v_U\circ w=j^\circ_A$. To this aim, it suffices to
show that $\bar v_U\circ w\circ u_{\bar A}=
j^\circ_A\circ u_{\bar A}:W(\bE(\bar A))\to A^\circ_U$.
However, by construction we have
$$
\bar v_U\circ w\circ u_{\bar A}=
v_U\circ W(j^\circ_\bE\circ\bar\omega^{-1})
\qquad\text{and}\qquad
j^\circ_A\circ u_{\bar A}\circ W(\bar\omega)=v_U\circ W(j^\circ_\bE)
$$
whence the claim.
\end{pfclaim}

(ii): Let $D\subset A_U^\circ$ be any open integrally
closed subring; it is easily seen that
$A_U^{\circ\circ}\subset D$, so $A_U^{\circ\circ}$ is an ideal
of $D$ (remark \ref{rem_someth-on-bdd-in-Z-lin}(iv))
and we set $\bar D:=D/A_U^{\circ\circ}\subset
C_A:=A_U^\circ/A_U^{\circ\circ}$.

\begin{claim}\label{cl_open-int-closed-subrings}
With the foregoing notation, the rule : $D\mapsto\bar D$
establishes a bijection from the set of open integrally
closed subrings of $A^\circ_U$ to the set of integrally
closed subrings of $C_A$.
\end{claim}
\begin{pfclaim} Let us check that $\bar D$ is integrally
closed in $C_A$. Indeed, let $\bar x\in C_A$, and
$P(T)\in D[T]$ a monic polynomial such that $P(\bar x)=0$
in $C_A$; if $x\in A_U^\circ$ is any representative for the
class $\bar x$, we get $P(x)\in A_U^{\circ\circ}$.
Set $Q(T):=P(T)-P(x)$; then $Q(T)\in D[T]$ and $Q(x)=0$,
so $x\in D$, and therefore $\bar x\in\bar D$, which shows
the contention. Conversely, if $\bar D\subset C_A$ is any
integrally closed subring, let us show that the preimage
$D\subset A_U^\circ$ of $\bar D$ is integrally closed
in $A_U^\circ$. Indeed, say that $x\in A_U^\circ$ is
integral over $D$; then the class $\bar x\in C_A$ of
$x$ is integral over $\bar D$, hence $\bar x\in\bar D$,
and the claim follows.
\end{pfclaim}

Claim \ref{cl_open-int-closed-subrings} also implies that
if $D\subset\bE_U^\circ$ is any open and integrally closed
subring, then $\bE_U^{\circ\circ}\subset D$, and the rule
$D\mapsto D/\bE_U^{\circ\circ}$ establishes a bijection
between the set of all open integrally closed subrings
of $\bE_U^\circ$ and the set of all integrally closed
subrings of $C_\bE:=\bE_U^\circ/\bE_U^{\circ\circ}$. However,
proposition \ref{prop_general-case}(iv) yields a natural
ring isomorphism $\bar\phi{}^{\flat\circ}_U:C_A\isom C_\bE$,
whence a bijection $\cS_A\isom\cS_\bE$, and it remains
to check that this bijection is realized by the rule
$D\mapsto\bE(D)$. To this aim, notice that every
integrally closed subring of $C_A$ is a perfect
$\F_p$-algebra; then the assertion follows easily
from proposition \ref{prop_construct-new-perfs},
which also gives (iii).
\end{proof}

\sset\subsubsection{}\label{subsec_renzi-victor}
In the situation of \eqref{subsec_general-case}, pick
a finitely generated ideal $J_\bE\subset\bE$ with
$\alpha_0\in J_\bE$, and such that $Z_\bE=\Spec\,\bE/J_\bE$.
In light of \eqref{eq_better}, it follows easily that
$Z_A=\Spec\,A/J_\bE^{\La 1\Ra}A$. Moreover, it is easily
seen that $J_\bE^{\lceil 0\rceil}$ is the radical of $J_\bE$;
on the other hand, the map $\bar u_A$ induces an
isomorphism $\bE/J_\bE^{\lceil 0\rceil}\isom A/J_\bE^{\lceil 0\rceil}A$,
so $A/J_\bE^{\lceil 0\rceil}A$ is reduced, {\em i.e.}
$J_\bE^{\lceil 0\rceil}A$ is a radical ideal of $A$.

\begin{corollary}\label{cor_renzi-victor}
With the notation of \eqref{subsec_renzi-victor},
suppose moreover that $A=A^+$. Then :
$$
J_\bE^{\lceil 0\rceil}A{}^+_U=J_\bE^{\lceil 0\rceil}\bar A.
$$
\end{corollary}
\begin{proof} Notice that it suffices to show that
$J_\bE^{\lceil 0\rceil}A{}^+_U\subset\bar A$. Indeed, suppose
that the latter holds; since we have as well
$(J_\bE^{\lceil 0\rceil})^2=J_\bE^{\lceil 0\rceil}$, we deduce that
$J_\bE^{\lceil 0\rceil}A{}^+_U\subset J_\bE^{\lceil 0\rceil}\bar A$,
and the converse inclusion is obvious. Notice that $\bE^+_U$
is the integral closure of $\bar\bE$ in $\bE_U$, so that
theorem \ref{th_perfect-purity}(ii) implies that the
assertion holds for $\bE^+_U$, {\em i.e.}
\set\begin{equation}\label{eq_hence}
J^{\lceil 0\rceil}_\bE\bE^+_U\subset\bar\bE.
\end{equation}
Next, recall that proposition \ref{prop_general-case}(iv)
yields a ring isomorphism
$\bar\phi_U^{\flat\circ}:\bE_U^\circ/\bE^{\circ\circ}\isom
A_U^\circ/A_U^{\circ\circ}$ restricting to an isomorphism
$\bar\bE/\bE_U^{\circ\circ}\isom\bar A/A_U^{\circ\circ}$.
Moreover, claim \ref{cl_open-int-closed-subrings} shows
that the subring $\bE^+_U/\bE_U^{\circ\circ}$ is the integral
closure of $\bar\bE/\bE_U^{\circ\circ}$ in
$\bE_U^\circ/\bE^{\circ\circ}$ and likewise,
$A^+_U/A_U^{\circ\circ}$ is the integral closure
of $\bar A/A_U^{\circ\circ}$ in $A_U/A_U^{\circ\circ}$. Hence,
$A^+_U/A_U^{\circ\circ}=
\bar\phi_U^{\flat\circ}(\bE^+_U/\bE_U^{\circ\circ})$, and taking
into account \eqref{eq_hence}, the assertion follows easily.
\end{proof}

\sset\subsubsection{}\label{subsec_Gamma-circ}
Keep the notation of \eqref{subsec_general-case}.
By virtue of lemma \ref{lem_dots-on-is}(i,iii) we may
identify naturally $U_A$ with an open subset of
$X_A^\circ:=\Spec\,A^\circ_U$, and we have a well defined
quasi-affinoid ring
$\underline A{}^\circ_U:=(A_U^\circ,A^+_U,U_A)$, which is
perfectoid, by theorem \ref{th_int-subrings-perfectoid}(i).
Set
$$
(U_A,\cT^\circ_{U_A},A^+_U):=\sSpec\,\underline A{}^\circ_U.
$$
since, by definition, the continuous maps $A\to A^\circ_U\to A_U$
are both f-adic, we see that the topology $\cT^\circ_{U_A}$ agrees
with $\cT_{U_A}$; {\em i.e.} we have a natural identification
$$
\sSpec\,\underline A{}^\circ_U\isom\sSpec\,\underline A.
$$
Lastly, it is clear that $\underline A{}^\circ_U$ depends
only on $\sSpec\,\underline A$, and we claim that the
rule $\sSpec\,\underline A\mapsto\underline A{}^\circ_U$
extends to a well defined functor
$$
\sGamma^\circ:\mathsf{q.Afd.Sch}_\mathrm{perf}\to
\mathsf{q.Afd.Ring}^o_\mathrm{perf}
$$
with a natural isomorphism $\eps_{\underline X}:
\sSpec\circ\sGamma^\circ(\underline X)\isom\underline X$
for every object
$\underline X$ of $\mathsf{q.Afd.Sch}_\mathrm{perf}$.
Indeed, if $\underline X:=(X,\cT_X,A^+_X)$ and
$\underline Y:=(Y,\cT_Y,A^+_Y)$ are any two perfectoid
quasi-affinoid schemes, and
$\psi:\underline X\to\underline Y$ any f-adic morphism
of quasi-affinoid schemes, then the corresponding map
$\psi^\flat:\cO_Y(Y)\to\cO_{\!X}(X)$ induces a morphism
$\psi^{\flat\circ}:(\cO_Y(Y)^\circ,A^+_Y)\to(\cO_{\!X}(X)^\circ,A^+_X)$
of affinoid rings (lemma \ref{lem_f-adics}(iii.a)), which is
obviously also f-adic, and the resulting diagram
$$
\xymatrix{
\Spec\,\cO_{\!X}(X) \ar[rr]^-{\Spec\,\psi^\flat} \ar[d]_{j_X} & &
\Spec\,\cO_Y(Y) \ar[d]^{j_Y} \\
\Spec\,\cO_{\!X}(X)^\circ \ar[rr]^-{\Spec\,\psi^{\flat\circ}} & &
\Spec\,\cO_Y(Y)^\circ
}$$
shows that $\Spec\,\psi^{\flat\circ}$ restricts to a morphism
$j_X(X)\to j_Y(Y)$, whence a well defined morphism
of perfectoid quasi-affinoid rings :
$$
\sGamma^\circ(\psi):
\sGamma^\circ(\underline Y)\to\sGamma^\circ(\underline X).
$$
Moreover, if $\underline A:=(A,A^+,U)$ is any perfectoid
quasi-affinoid ring, the natural map
$A\to\cO_{\!U}(U)^\circ$ induces a well defined morphism
of quasi-affinoid rings $\eta_{\underline A}:
\underline A\to\sGamma^\circ\circ\sSpec(\underline A)$
which is f-adic, and therefore it is a morphism of
perfectoid quasi-affinoid rings. Lastly, it is easily
seen that the pair $(\eps_\bullet,\eta_\bullet)$ fulfills
the triangular conditions of \eqref{subsec_adj-pair},
so $\sGamma^\circ$ is indeed right adjoint to $\sSpec$.

\sset\subsubsection{}\label{subsec_upgrade-E}
We also want to upgrade the functor $\bE$ to well
defined functors
$$
\bE:\mathsf{q.Afd.Ring}_\mathrm{perf}\to
\mathsf{q.Afd.Ring}_\mathrm{perf}
\qquad
\bE:\mathsf{q.Afd.Sch}_\mathrm{perf}\to
\mathsf{q.Afd.Sch}_\mathrm{perf}.
$$
Namely, let $\underline A:=(A,A^+,U)$ be any perfectoid
quasi-affinoid ring, and set $\bE:=\bE(A)$; by applying
theorem \ref{th_int-subrings-perfectoid}(ii) to the
perfectoid quasi-affinoid ring $(A,A^+,\Spec\,A)$, we see
that $\bE^+:=\bE(A^+)$ is a subring of integral elements
of $\bE$, so we get a perfectoid quasi-affinoid ring
$$
\bE(\underline A):=(\bE,\bE^+,\bE(U))
$$
where $\bE(U)\subset\Spec\,\bE$ is defined as in
\eqref{subsec_general-case} : namely, it is the unique
open subset of $\Spec\,\bE$ containing the analytic locus
and such that $\phi^{-1}\bE(U)=U$, where
$\phi:\Spec\,A\to\Spec\,\bE$ is the continuous map given by
\eqref{subsec_map-on-spec}. Moreover, if $\underline\alpha$
is any distinguished element in $\Ker\,u_A$, endow $A_0:=A/pA$
and $\bE_0:=\bE/\alpha_0\bE$ with the ring topologies induced
by $A$ and respectively $\bE$; from remark
\ref{rem_nice-topology}(ii) we see as well that the map
$u_A:A\to\bE$ induces a natural isomorphism of quasi-affinoid
rings (notation of example \ref{ex_shouldbe-quasi-affinoid}(i)) :
\set\begin{equation}\label{eq_upgrade-rem-nice-ii}
A_0\otimes_A\underline A\isom\bE_0\otimes_\bE\bE(\underline A).
\end{equation}
Furthermore, for every morphism
$f:\underline A\to\underline B:=(B,B^+,V)$ of perfectoid
quasi-affinoid rings, it follows easily from
\eqref{eq_functoriality-of-Spec_u} that the morphism of
schemes $\Spec\,\bE(f):\Spec\,\bE(B)\to\Spec\,\bE$ maps
$\bE(V)$ into $\bE(U)$. Moreover, the ring homomorphism
$\bE(f):\bE\to\bE(B)$ is adic (theorem \ref{th_adic-to-adic}(i)),
whence a well defined morphism of perfectoid quasi-affinoid rings
$$
\bE(f):\bE(\underline A)\to\bE(\underline B).
$$
Let also $\underline B\tdu\mathsf{q.Afd.Ring}_\mathrm{perf}:=
\underline B/\mathsf{q.Afd.Ring}_\mathrm{perf}$ for every
perfectoid quasi-affinoid ring $\underline B$; just as in
remark \ref{rem_nice-topology}(i), this new functor $\bE$
restricts to an equivalence
$$
\underline A\tdu\mathsf{q.Afd.Ring}_\mathrm{perf}\isom
\bE(\underline A)\tdu\mathsf{q.Afd.Ring}_\mathrm{perf}
$$
which admits a natural quasi-inverse functor
$$
\bA:\bE(\underline A)\tdu\mathsf{q.Afd.Ring}_\mathrm{perf}
\isom\underline A\tdu\mathsf{q.Afd.Ring}_\mathrm{perf}
\qquad
(E,E^+,V)\mapsto(\bA(E),\bA(E^+),\bA(V))
$$
with $\bA(E):=W(E)\otimes_{W(A)}A$ and correspondingly for
$\bA(E^+)$, and where $\bA(V):=\phi^{-1}V$. Likewise, for
every perfectoid quasi-affinoid scheme $\underline U$ we
define
$$
\bE(\underline U):=
\sSpec\,\circ\bE\circ\sGamma^\circ(\underline U)
\qquad\text{and}\qquad
\underline U\tdu\mathsf{q.Afd.Sch}_\mathrm{perf}:=
\mathsf{q.Afd.Sch}_\mathrm{perf}/\underline U.
$$
Then, the resulting functor $\bE$ on quasi-affinoid schemes
restricts to an equivalence
$$
\underline U\tdu\mathsf{q.Afd.Sch}_\mathrm{perf}\isom
\bE(\underline U)\tdu\mathsf{q.Afd.Sch}_\mathrm{perf}
$$
with quasi-inverse given by the rule :
$\underline V\mapsto\sSpec\circ\bA\circ\sGamma^\circ(\underline V)$
(details left to the reader). Lastly, say that
$\underline U=\sGamma^\circ(\underline A)$, and set
$\underline U_0:=\sSpec\,(A_0\otimes_A\underline A)$ and
$\bE(\underline U)_0:=\sSpec\,(\bE_0\otimes_\bE\underline\bE)$;
then we deduce from \eqref{eq_upgrade-rem-nice-ii} a natural
isomorphism of quasi-affinoid schemes
\set\begin{equation}\label{eq_osiris}
\underline U_0\times_{\underline U}\underline V\isom
\bE(\underline U)_0\times_{\bE(\underline U)}\bE(\underline V)
\qquad
\text{for every $\underline V\in
\Ob(\underline U\tdu\mathsf{q.Afd.Sch}_\mathrm{perf})$}.
\end{equation}
(see example \ref{ex_fibre-prod-in-qaff-sch}).

\sset\subsubsection{}\label{eq_back-on-track}
Next, we wish to extend the results of section
\ref{subsec_Hom-theory-perfectoid} to the perfectoid
quasi-affinoid case. Resume the notation of
\eqref{subsec_general-case}, and pick any finitely
generated ideal of definition $I$ of $\bE$; we get a
descending filtration on $\bE_U:=\cO_{\!U_\bE}(U_\bE)$
and $A_U:=\cO_{U_A}(U_A)$ by the rules :
$$
\Fil^s\bE_U:=I^{\La s\Ra}\bar\bE
\qquad
\Fil^sA_U:=I^{\La s\Ra}\bar A
$$
$$
\Fil^{-s}\bE_U:=\{x\in\bE_U~|~I^{\La s\Ra}x\subset\bar\bE\}
\qquad
\Fil^{-s}A_U:=\{x\in A_U~|~I^{\La s\Ra}x\subset\bar A\}
$$
for every $s\in\N[1/p]$. As in example \ref{ex_Samuel}(i)
we associate with $I$ an {\em angular order function}
$$
\nu:\bE_U\to\R\cup\{+\infty\}
\qquad
x\mapsto\sup\{s\in\Z[1/p]~|~x\in\Fil^s\bE_U\}.
$$
Notice that, since the topology of $\bE_U$ is separated,
we have $\nu(x)=+\infty$ if and only if $x=0$. Then, we fix
$\rho\in]0,1[$, and set
$$
|x|_I:=\rho^{\nu(x)}
\qquad
\text{for every $x\in\bE_U$}
$$
(with the convention that $\rho^{+\infty}:=0$).

\begin{lemma}\label{lem_Samuel-perfectoid}
With the notation of \eqref{eq_back-on-track}, we have :
\begin{enumerate}
\item
The mapping $|\cdot|_I$ is the asymptotic Samuel function
on $\bE_U$ associated with $I$ (see example
{\em\ref{ex_Samuel}(i)}). Especially, it is a real-valued
power-multiplicative norm on $\bE_U$.
\item
$\phi^\flat_U(\Fil^s\bE_U)\subset\Fil^sA_U$ for every $s\in\Z[1/p]$
(notation of proposition {\em \ref{prop_new-formula}(i)}).
\end{enumerate}
\end{lemma}
\begin{proof}(i): Let $|\cdot|^*_I$ be the asymptotic Samuel
function attached to $I$ (and the real number $\rho$), let
$x\in\bE_U$ be any element, and $s\in\Z[1/p]$ any rational
number; consider the following conditions:
\begin{enumerate}
\alphaenu
\item
$|x|_I\leq\rho^s$.
\item
For every $t<s$ in $\Z[1/p]$ we have $x\in\Fil^t\bE_U$.
\item
For every $t<s$ in $\Z[1/p]$ we may find $n\in\N$ such that
$p^nt\in\Z$ and $x^{p^n}\in I^{p^nt}$, where the integral power
$I^{p^nt}$ is defined as in example \ref{ex_Samuel}(i).
\item
$|x|^*_I\leq\rho^s$.
\end{enumerate}
It is easily seen that each of these conditions is equivalent
to the following one, whence the assertion.

(ii): Suppose first that $s\geq 0$; in this case, by
definition we have $\phi^\flat_U(I^{\La s\Ra})\subset\Fil^sA_U$.
On the other hand, proposition \ref{prop_new-formula}(ii)
implies that $\phi^\flat_U(\Fil^s\bE_U)$ lies in the topological
closure of the $\bar A$-submodule of $A_U$ generated by
$\phi^\flat_U(I^{\La s\Ra})$, and since $\Fil^sA_U$ is an open
$\bar A$-submodule of $A_U$, the assertion follows.

Lastly, suppose that $s<0$, and let $x\in\Fil^s\bE_U$ be any
element; the condition means that $I^{\La s\Ra}x\subset\bar\bE$,
whence $\phi^\flat_U(x)\cdot\phi_U^\flat(I^{\La s\Ra})\subset\bar A$,
and finally $\phi^\flat_U(x)\in\Fil^sA_U$, as stated.
\end{proof}

\sset\subsubsection{}\label{subsec_boeuf}
The norm $|\cdot|_I$ in turns induces a mapping
$$
|\cdot|_1:W(\bE_U)\to\R_+\cup\{+\infty\}
\qquad
(a_n~|~n\in\N)\mapsto\sup_{n\in\N}|a_n|_I^{p^{-n}}
$$
as in \eqref{subsec_seminorm-on-W}, such that the subset
$$
W(\bE_U,1):=\{\underline a\in W(\bE_U)~|~|\underline a|_1\in\R_+\}
$$
is a subring of $W(\bE_U)$, and the restriction of $|\cdot|_1$
is a real-valued norm on $W(\bE_U,1)$ (proposition
\ref{prop_semi-norm-on-W}(ii)). We endow $W(\bE_U,1)$
with the topology defined by the norm $|\cdot|_1$. Then
clearly $W(\bE_U,1)$ is a separated topological ring,
and it is independent of the choice of $I$ (see example
\ref{ex_Samuel}(iv)). Furthermore, for every $s\in\Z[1/p]$
we set
$$
W\La s\Ra:=\{(x_n~|~n\in\N)\in W(\bE_U)~|~
x_n\in\Fil^{p^ns}\bE_U\ \text{for every $n\in\N$}\}
\subset W(\bE_U,1).
$$
Recall that the f-adic topologies of $\bE_U$ and $A_U$
agree with the $\Z$-linear topologies defined by
$\Fil^\bullet\bE_U$ and respectively $\Fil^\bullet A_U$
(lemma \ref{lem_mon-fract-powers}(iv) and corollary
\ref{cor_taut-two}(ii)). Moreover, $\Fil^sA_U$ and
$\Fil^s\bE_U$ are bounded subsets of $A_U$ and respectively
$\bE_U$, for every $s\in\Z[1/p]$. It follows easily that
for every $\underline a:=(a_n~|~n\in\N)\in W(\bE_U,1)$,
the series
$$
\sum_{n\in\N}p^n\cdot\phi^\flat_U(a_n^{1/p^n})
$$
converges to a well defined element $u_U(\underline a)$ of
$A_U$. Moreover, just as in definition \ref{def_beta-taut}(v),
for any $\bar\bE$-submodule $\cK$ of $\bE_U$ we denote by
$\{\cK\}$ the topological closure in $A_U$ of the
$\bar A$-submodule generated by $(\phi^\flat_U(x)~|~x\in\cK)$.
We also define filtrations on $\cK$ and $\{\cK\}$ by the
rules :
$$
\Fil^s\cK:=\cK\cap\Fil^s\bE_U
\qquad\text{and}\qquad
\Fil^s\{\cK\}:=\{\Fil^s\cK\}
\qquad
\text{for every $s\in\Z[1/p]$}.
$$

\begin{proposition}\label{prop_boeuf}
With the notation of \eqref{subsec_boeuf}, we have :
\begin{enumerate}
\item
The resulting mapping $u_U:W(\bE_U,1)\to A_U$
is a morphism of topological rings.
\item
$u_U(W\La s\Ra)=\{\Fil^s\bE_U\}\subset\Fil^sA_U$ for every
$s\in\Z[1/p]$.
\item
$pW(\bE_U)\cap W\La s\Ra=pW\La s\Ra$ for every $s\in\Z[1/p]$.
\item
The topological ring $W(\bE_U,1)$ is complete and separated.
\end{enumerate}
\end{proposition}
\begin{proof}(i): The continuity of $u_U$ follows from the
continuity of $\phi_U^\flat$ (proposition
\ref{prop_new-formula}(i)) and a simple inspection of
the definition. Next, we remark :

\begin{claim}\label{cl_alone-in-here}
For every $a\in\bE_U$ and every $\underline x\in W(\bE_U,1)$
we have
$$
\tau_{\bE_U}(a)\cdot\underline x\in W(\bE_U,1)
\qquad\text{and}\qquad
u_U(\tau_{\bE_U}(a)\cdot\underline x)=
\phi^\flat_U(a)\cdot u_U(\underline x).
$$
\end{claim}
\begin{pfclaim} Say that $\underline x=(x_n~|~n\in\N)$;
by proposition \ref{prop_Teich-series}(i) we have
$\tau_{\bE_U}(a)\cdot\underline x=(a^{p^n}x_n~|~n\in\N)$,
whence the first assertion. By the same token, since
$\phi^\flat_U$ is a continuous morphism of multiplicative
monoids, we may compute
$$
u_U(\tau_{\bE_U}(a)\cdot\underline x)=
\sum_{n\in\N}p^n\cdot\phi^\flat_U(ax_n^{1/p^n})=
\phi_U^\flat(a)\cdot
\sum_{n\in\N}p^n\cdot\phi^\flat_U(x_n^{1/p^n})=
\phi^\flat_U(a)\cdot u_U(\underline x)
$$
as stated.
\end{pfclaim}

To conclude, we have to check that
$$
u_U(\underline x+\underline y)=
u_U(\underline x)+u_U(\underline y)
\qquad\text{and}\qquad
u_U(\underline x\cdot\underline y)=
u_U(\underline x)\cdot u_U(\underline y)
$$
for every $\underline x:=(x_n~|~n\in\N),
\underline y:=(y_n~|~n\in\N)\in W(\bE_U,1)$.
However, pick $s\in\Z[1/p]$ such that
$\underline x,\underline y\in W\La s\Ra$, and notice
that $u_U$ restricts to a map $W\La s\Ra\to\Fil^sA_U$; set
$\underline x^{(k)}:=\sum_{n=0}^kp^n\cdot\tau_{\bE_U}(x_n)$
and define likewise $\underline y^{(k)}$ for every $k\in\N$.
The sequences $(\underline x^{(k)}~|~n\in\N)$ and
$(\underline y^{(k)}~|~n\in\N)$ lie in $W\La s\Ra$
and converge $p$-adically to $\underline x$ and
respectively $\underline y$; moreover, the $p$-adic
topology is separated on $\Fil^sA_U$, so it suffices to check
the sought identities with $\underline x$ and $\underline y$
replaced by $\underline x^{(k)}$ and $\underline y^{(k)}$,
for every $k\in\N$. We may thus assume that there exists
$k\in\N$ such that $x_n=y_n=0$ for every $n>k$. Then, we
argue as in the proof of proposition \ref{prop_new-formula}(ii) :
we pick a family $(g_\lambda~|~\lambda\in\Lambda)$ of elements
of $\bE$ such that
$U_\bE=\bigcup_{\lambda\in\Lambda}\Spec\,\bE[g^{-1}_\lambda]$,
and we set
$$
h_\lambda:=\bar u_A(g_\lambda)\in A
\qquad\text{and}\qquad
t_\lambda:=\tau_\bE(g_\lambda)\in W(\bE)
\qquad
\text{for every $\lambda\in\Lambda$}
$$
so that $U_A=\bigcup_{\lambda\in\Lambda}\Spec\,A[h^{-1}_\lambda]$,
and it suffices to check that the sought identities hold
in $A[h_\lambda^{-1}]=\bar A[h^{-1}_\lambda]$, for every
$\lambda\in\Lambda$. Now, for every $\lambda\in\Lambda$
we may find $n_\lambda\in\N$ such that
$g_\lambda^{p^in_\lambda}\cdot x_i,g_\lambda^{p^in_\lambda}\cdot y_i
\in\bar\bE$ for every $i=0,\dots,k$, and therefore
$$
\underline x{}_\lambda:=t_\lambda^{n_\lambda}\cdot\underline x,\ 
\underline y{}_\lambda:=t^{n_\lambda}_\lambda\cdot\underline y
\in W(\bar\bE)
$$
(proposition \ref{prop_Teich-series}(i)). Lastly, we
consider the commutative diagram
$$
\xymatrix{ W(\bar\bE) \ar[r]^-{u_{\bar A}} \ar[d] &
\bar A \ar[r] \ar[d]
& \bar A[h^{-1}_\lambda] \ddouble \\
W(\bE_U,1) \ar[r]^-{u_U} & A_U \ar[r] & A[h^{-1}_\lambda]
}$$
whose two vertical arrows are the natural inclusions,
and whose unmarked horizontal arrows are the localization
maps; taking into account claim \ref{cl_alone-in-here},
we may compute for every $\lambda\in\Lambda$:
$$
\begin{aligned}
\phi^\flat_U(t_\lambda^{n_\lambda})\cdot u_U(\underline x+\underline y)
=&\, u_U(\underline x{}_\lambda+\underline y{}_\lambda) \\
=&\, u_{\bar A}(\underline x{}_\lambda+\underline y{}_\lambda) \\
=&\, u_{\bar A}(\underline x{}_\lambda)+
u_{\bar A}(\underline y{}_\lambda) \\
=&\, u_U(\underline x{}_\lambda)+u_U(\underline y{}_\lambda) \\
=&\, \phi^\flat_U(t_\lambda^{n_\lambda})\cdot
(u_U(\underline x)+u_U(\underline y))
\end{aligned}
$$
as required. Likewise one shows the other sought identity.

(ii): Notice first that if $s\geq 0$, then
$W\La s\Ra$ is the ideal $W(I^{\La s\Ra}\bar\bE)$ of
$W(\bar\bE)$ (notation of remark \ref{rem_Witt-limit}(iv));
since $I^{\La s\Ra}\bar\bE$ is a taut and open ideal of
$\bar\bE$, the ideal $W(I^{\La s\Ra}\bar\bE)$ is closed
in $W(\bar\bE)$, and from theorem \ref{th_taut-two}(i)
it follows that
$$
u_U(W\La s\Ra)=u_{\bar A}(W(I^{\La s\Ra}\bar\bE))=\{\Fil^s\bE_U\}
$$
as required. Especially $u_U(W\La s\Ra)$ is an open
$\bar A$-submodule of $A_U$ for every $s\geq 0$;
but then obviously the same holds also more generally
for every $s\in\Z[1/p]$. Next, for a general $s\in\Z[1/p]$
notice that every element of $W\La s\Ra$ can be written
as a $p$-adically convergent series
$\sum_{n\in\N}p^n\cdot\tau_{\bE_U}(x_n)$, for a unique sequence
$(x_n~|~n\in\N)$ of elements of $\Fil^s\bE_U$; to any such
element, the map $u_U$ assigns the $p$-adically convergent
series $\sum_{n\in\N}p^n\cdot\phi^\flat_U(x_n)$, which clearly
lies in $\{\Fil^s\bE_U\}$, from which we see that
$u_U(W\La s\Ra)$ is dense in $\{\Fil^s\bE_U\}$. However,
we have just seen that $u_U(W\La s\Ra)$ is a closed subset
of $A_U$ for every such $s$, whence the assertion.

(iii): The intersection $pW(\bE_U)\cap W\La s\Ra$
consists of all elements $\underline x:=(x_n~|~n\in\N)$
such that $x_0=0$ and $x_n\in\Fil^{sp^n}\bE_U$ for every
$n\geq 1$. For such $\underline x$, we have
$y_n:=x_{n+1}^{1/p}\in\Fil^{sp^n}\bE_U$ for every $n\in\N$,
whence $\underline y:=(y_n~|~n\in\N)\in W\La s\Ra$ and
$p\cdot\underline y=\underline x$, whence the contention.

(iv): Since $W(\bar\bE)$ is an open subring of $W(\bE_U,1)$,
it suffices to show that $W(\bar\bE)$ is complete and
separated for the topology $\cT$ induced by the norm
$|\cdot|_1$. To this aim, let $x_1,\dots,x_n$ be any
finite system of generators of $I\bar\bE$, and denote
by $\cI\subset W(\bar\bE)$ the ideal generated by
$\tau_{\bar\bE}(x_1),\dots,\tau_{\bar\bE}(x_n)$. Taking into
account proposition \ref{prop_morel}(i) and lemma
\ref{lem_mon-fract-powers}(iv), it is easily seen that
$\cT$ agrees with the $\cJ$-adic topology on $W(\bar\bE)$.
However, let also $\cT_{\bar\bE}$ be the $I\bar\bE$-adic
topology on $\bar\bE$ and $\cT_{W(\bar\bE)}$ the topology of
$W(\bar\bE,\cT_E)$ (as in definition \ref{def_Witt-vectors});
then $\cT_{W(\bar\bE)}$ agrees with the $(pW(\bar\bE)+\cI)$-adic
topology (proposition \ref{prop_morel}(ii)). On the other
hand, since $\cT_{\bar\bE}$ is complete and separated, the same
holds for $\cT_{W(\bar\bE)}$ (lemma \ref{lem_Witt-limit}(ii)).
To conclude, we may now appeal to lemma \ref{lem_fontaine}.
\end{proof}

\sset\subsubsection{}\label{subsec_quasi-taut}
In the situation of \eqref{subsec_boeuf}, let now $\beta\in\bE$
be any element, and $\cJ\subset\cK$ two $\bar\bE$-submodules of
$\bE_U$; just as in definition \ref{def_beta-taut}, we shall say
that the inclusion of $\cJ$ in $\cK$ is {\em $\beta$-taut} if
$$
\beta\cdot\Phi_{\bE_U}^{-1}(\cK^p)\subset\cJ
$$
where $\Phi_{\bE_U}$ is the Frobenius automorphism of $\bE_U$.
For $\beta=\alpha_0$ (where $(\alpha_n~|~n\in\N)$ is a fixed
distinguished element in $\Ker\,u_A$), we just say that the
inclusion is {\em taut}. Likewise, we say that $\cK$ is
{\em $\beta$-taut} (resp. {\em taut}) if the identity map of
$\cK$ is a $\beta$-taut (resp. taut) inclusion.
Just as in \eqref{subsec_elementary}, we see that for any
taut inclusion $\cJ\subset\cK$, the quotients $\cK/\cJ$ and
$\{\cK\}/\{\cJ\}$ are both $\bar A/p\bar A$-modules. Also,
just as in remark \ref{rem_beta-taut}(v), if the inclusion
$\cJ\subset\cK$ is taut, then both $\cJ$ and $\cK$ are taut.
With this terminology, we have the following extension of
theorem \ref{th_taut-one} :

\begin{lemma}\label{lem_boody-boody}
In the situation of \eqref{subsec_boeuf}, the following holds :
\begin{enumerate}
\item
Every taut inclusion $\cK_1\subset\cK_2$ of\/
$\bar\bE$-submodules of\/ $\bE_U$ induces an
$\bar A/p\bar A$-linear map
$$
\cK_2/\cK_1\to\{\cK_2\}/\{\cK_1\}
\quad :\quad
(x\mod{\cK_1})\mapsto(\phi^\flat_U(x)\mod{\cK_1}).
$$
\item
If $\cK_1$ and $\cK_2$ are topologically closed
in $\bE_U$, the map of\/ {\em (i)} is an isomorphism.
\item
For every taut $\bar\bE$-submodule $\cK$ of\/ $\bE_U$ and
every $t\in\Z[1/p]$ we have
$$
\{\cK\}\cap\{\Fil^t\bE_U\}=\Fil^t\{\cK\}.
$$
\end{enumerate}
\end{lemma}
\begin{proof}(i): The proof is the same as that of theorem
\ref{th_taut-one}(i) : the only difference is that instead
of using proposition \ref{prop_combinatorial} one must
appeal to the more general proposition
\ref{prop_new-formula}(ii).

Next, we prove the following special case of (ii) :

\begin{claim}\label{cl_spec-case-qaff}
Pick any strictly positive $\eps\in\N[1/p]$ such that
$\alpha_0\in I^{\La\eps\Ra}\bE$. Then we have :
\begin{enumerate}
\item
For every $s,t\in\Z[1/p]$ such that $t-\eps\leq s\leq t$, the
inclusion $\Fil^t\bE_U\subset\Fil^s\bE_U$ is taut, and the map
of (i) is an isomorphism
$$
\tau_{s,t}:
\Fil^s\bE_U/\Fil^t\bE_U\isom\{\Fil^s\bE_U\}/\{\Fil^t\bE_U\}.
$$
\item
Moreover, the inclusion $\{\Fil^s\bE_U\}\subset\Fil^sA_U$
induces an injection
$$
\mu_{s,t}:\{\Fil^s\bE_U\}/\{\Fil^t\bE_U\}\to\Fil^sA_U/\Fil^tA_U.
$$
\end{enumerate}
\end{claim}
\begin{pfclaim}(i): First notice that, since $\bar A$ is
perfectoid and $\bE(\bar A)=\bar\bE$, the assertion in case
$s\geq 0$ is a special case of theorem \ref{th_taut-one}(ii).
Thus, we assume henceforth that $s<0$. Next, suppose that
$t>0$; in this case, we know already that $\tau_{0,t}$ is an
isomorphism, and therefore the $5$-lemma implies that
$\tau_{s,t}$ is an isomorphism if and only if the same
holds for $\tau_{s,0}$ (details left to the reader).
Hence, we may further assume that $t\leq 0$ as well.
Then, notice that the map $\tau_{s,t}$ has dense image,
for the quotient topology on the target; however,
$\{\Fil^t\bE_U\}$ is an open $\bE_U$-submodule of
$\{\Fil^t\bE_U\}$, so this quotient topology is discrete,
and thus $\tau_{s,t}$ is surjective. To conclude, it then
suffices to check that $\mu_{s,t}\circ\tau_{s,t}$ is injective.
Hence, let $x\in\Fil^s\bE_U$ be any element, and suppose
that $\phi^\flat_U(x)\in\Fil^tA_U$; we need to show that
$x\in\Fil^t\bE_U$, {\em i.e.} that for every $b\in I^{\La-t\Ra}$
we have $z:=bx\in\bar\bE$. However, for such $z$ we have
$$
\phi^\flat_U(z)=\phi^\flat_U(b)\cdot\phi_U^\flat(x)\in
I^{\La-t\Ra}\cdot\Fil^tA_U\subset\bar A.
$$
Notice that $z\in\Fil^{s-t}\bE_U$, and let
$\bar z\in\Fil^{s-t}\bE_U/\Fil^0\bE_U=\Fil^{s-t}\bE_U/\bar\bE$
be the class of $z$. Since $s-t\geq -\eps$, and since
$\alpha_0\in I^{\La\eps\Ra}\bE$, we see that
$\Fil^{s-t}\bE_U/\bar\bE\subset\Ann_{\bE_U/\bar\bE}(\alpha_0)$;
on the other hand, proposition \ref{prop_Q-c-tricky}(iv)
implies that $\phi^\flat_U$ induces an isomorphism
$$
\Ann_{\bE_U/\bar\bE}(\alpha_0)\isom\Ann_{A_U/\bar A}(p)
$$
whence the contention. Assertion (ii) is an immediate
consequence.
\end{pfclaim}

(iii): Let $\eps,s,t\in\Z[1/p]$ be as in claim
\ref{cl_spec-case-qaff}(i), and notice that the composition
$$
M:=\frac{\Fil^s\cK}{\Fil^t\cK}\xrightarrow{\ \beta\ }
M':=\frac{\Fil^s\{\cK\}}{\Fil^t\{\cK\}}
\xrightarrow{\ \gamma\ }
\frac{\{\cK\}\cap\{\Fil^s\bE_U\}}{\{\cK\}\cap\{\Fil^t\bE_U\}}\to
M'':=\frac{\{\Fil^s\bE_U\}}{\{\Fil^t\bE_U\}}
$$
factors through an injective map $M\to\Fil^s\bE_U/\Fil^t\bE_U$
and the isomorphism $\tau_{s,t}$. Hence, $\beta$ is injective;
moreover, $\beta$ has dense image, and the quotient topologies
on $M$ and $M''$ are discrete, so $\beta$ is an isomorphism,
by claim \ref{cl_little-top-trick}. We deduce that $\gamma$
is injective.

\begin{claim}\label{cl_ok-if-bounded}
Assertion (iii) holds if there exists $s\in\Z[1/p]$
such that $\cK\subset\Fil^s\bE_U$.
\end{claim}
\begin{pfclaim} Notice that $\cK\cap\Fil^r\bE_U$ is still
taut for every $r\in\Z[1/p]$. Now, if $t<s$ there is nothing
to prove. In case $t\geq s$, let $n\in\N$ be the unique integer
such that $t-n\eps-s\in[0,\eps[$; arguing by induction on $n$,
we are easily reduced to the case where $n=0$, {\em i.e.} we
may assume that $t-\eps<s\leq t$. In this case, clearly $\gamma$
is an isomorphism, and the sought identity follows immediately.
\end{pfclaim}

Now, for a general $\cK$, clearly $\{\cK\}$ is the topological
closure of $\bigcup_{s\in\Z[1/p]}\{\cK\cap\Fil^s\bE_U\}$ in $A_U$,
and since $\{\Fil^t\bE_U\}$ is open in $A_U$, we deduce that
$\{\cK\}\cap\{\Fil^t\bE_U\}$ is the topological closure of
$\bigcup_{s\in\Z[1/p]}(\{\cK\cap\Fil^s\bE_U\}\cap\{\Fil^t\bE_U\})$.
But for every $s\leq t$, we have
$\{\cK\cap\Fil^s\bE_U\}\cap\{\Fil^t\bE_U\}=\{\cK\cap\Fil^t\bE_U\}$,
by virtue of claim \ref{cl_ok-if-bounded}, whence (iii).

We may now complete the proof of (ii) arguing as in the
proof of theorem \ref{th_taut-one}(ii) : namely, we may
assume that $\alpha_0\in I$, and for every $\bar\bE$-submodule
$\cK$ of $\bE_U$ we set
$$
\gr^n\cK:=\Fil^n\cK/\Fil^{n+1}\cK
\qquad
\gr^n\{\cK\}:=\Fil^n\{\cK\}/\Fil^{n+1}\{\cK\}
\qquad
\text{for every $n\in\Z$}.
$$
Moreover, for every $n\in\Z$ we let $\Fil^n(\cK_2/\cK_1)$
(resp. $\Fil^n(\{\cK_2\}/\{\cK_1\})$) be the image of
$\Fil^n\cK_2$ in $\cK_2/\cK_1$ (resp. of
$\Fil^b\{\cK_2\}$ in $\{\cK_2\}/\{\cK_1\}$) and we denote
again by $\gr^\bullet(\cK_2/\cK_1)$ and
$\gr^\bullet(\{\cK_2\}/\{\cK_1\})$ the respective associated
graded modules. There follows for every $n\in\Z$ a commutative
diagram with exact rows
$$
\xymatrix{
0 \ar[r] & \gr^n\cK_1 \ar[r] \ar[d] &
\gr^n\cK_2 \ar[r] \ar[d] &
\gr^n(\cK_2/\cK_1) \ar[r] \ar[d] & 0 \\
0 \ar[r] & \gr^n\{\cK_1\} \ar[r] & \gr^n\{\cK_2\} \ar[r] &
\gr^n(\{\cK_2\}/\{\cK_1\}) \ar[r] & 0
}$$
and since both $\cK_1$ and $\cK_2$ are taut, we know
already from the proof of (iii) that the leftmost and
central vertical arrows are isomorphisms, hence the
same holds for the rightmost vertical arrow. Now, notice
that the filtration $\Fil^\bullet(\cK_2/\cK_1)$ defines a
separated and complete topology on $\cK_2/\cK_1$, since
$\cK_2$ and $\cK_1$ are both closed $\bar\bE$-submodules
of $\bE$. The same holds for the topology on
$\{\cK_2\}/\{\cK_1\}$ determined by the filtration
$\Fil^\bullet(\{\cK_2\}/\{\cK_1\})$, since $\{\cK_2\}$
and $\{\cK_1\}$ are closed $\bar A$-submodules in $A$.
Then (ii) follows directly from
\cite[Ch.III, \S2, n.8, Cor.3]{BouAC}.
\end{proof}

\begin{theorem}\label{th_boody-boody}
In the situation of \eqref{subsec_boeuf}, we have :
\begin{enumerate}
\item
The continuous ring homomorphism $u_U$ is surjective and open.
\item
$\Ker\,u_U=\underline\alpha W(\bE_U,1)$.
\end{enumerate}
\end{theorem}
\begin{proof}(ii): Clearly $\underline\alpha\in\Ker u_U$;
thus, it suffices to check that $u_U$ induces an isomorphism
$$
\bar u_s:
W\La s\Ra/\underline\alpha W\La s\Ra\isom\{\Fil^s\bE_U\}
\qquad
\text{for every $s\in\Z[1/p]$}
$$
and by virtue of proposition \ref{prop_boeuf}(ii) we know
already that $\bar u_s$ is surjective for every such $s$.
Now, if $s\geq 0$, we have
$W\La s\Ra=W(I^{\La s\Ra}\bar\bE)\subset\W(\bar\bE)$, and the
assertion follows from theorem \ref{th_taut-two}(iv).
For the general case, we may assume that $\alpha_0\in I$,
and then an easy induction argument then reduces to checking :

\begin{claim} Suppose that $\alpha_0\in I$. Then, if
$\bar u_{s+1}$ is injective, the same holds for $\bar u_s$.
\end{claim}
\begin{pfclaim} We consider the commutative diagram
\set\begin{equation}\label{eq_boody-boody}
{\diagram W\La s+1\Ra/\underline\alpha W\La s+1\Ra
\ar[r]^-j \ar[d]_{\bar u_{s+1}} &
W\La s\Ra/\underline\alpha W\La s\Ra \ar[d]^{\bar u_s} \\
\{\Fil^{s+1}\bE_U\} \ar[r]^-i & \{\Fil^s\bE_U\}
\enddiagram}
\end{equation}
where $i$ is the inclusion map, and $j$ is induced
by the inclusion $W\La s+1\Ra\subset W\La s\Ra$. Since
$\bar u_{s+1}$ is injective by assumption, the same holds
for $j$. Next, let us write
$\alpha=\tau_\bE(\alpha_0)+p\cdot\underline u$ for some
invertible element $\underline u$ of $W(\bE)$; from
proposition \ref{prop_Teich-series}(ii) and our
assumption on $\alpha_0$, we deduce that
$\tau_\bE(\alpha_0)\cdot W\La s\Ra\subset W\La s+1\Ra$,
hence the cokernel of $j$ is annihilated by $p$. Taking
into account proposition \ref{prop_boeuf}(iii), it follows
that $\Coker\,j$ is naturally identified with the cokernel
of the natural map
$j':(\Fil^{s+1}\bE_U)/\alpha_0(\Fil^{s+1}\bE_U)\to
(\Fil^s\bE_U)/\alpha_0(\Fil^s\bE_U)$.
But by the same token we see that
$\alpha_0(\Fil^s\bE_U)\subset\Fil^{s+1}\bE_U$, so
$\Coker\,j'=\Fil^s\bE_U/\Fil^{s+1}\bE_U$. Moreover,
a simple inspection of the definitions shows that
the map $\Coker\,j'\to\Coker\,i$ resulting from
\eqref{eq_boody-boody} is the same as the isomorphism
$\tau_{s,s+1}$ of claim \ref{cl_spec-case-qaff}(i).
Then the assertion follows from the $5$-lemma.
\end{pfclaim}

(i): By proposition \ref{prop_boeuf}(ii) we know already
that $u_U$ is an open map. Next, we remark :

\begin{claim} To prove that $u_U$ is surjective, we
may assume that the topology of $\bE$ is $\alpha_0$-adic.
\end{claim}
\begin{pfclaim} Notice that the image $R$ of $u_U$ contains
$u_U(W(\bar\bE))=u_{\bar A}(W(\bar\bE))=\bar A$ (lemma
\ref{lem_was-third-cond}(i)). Especially, $R$ is an
open subring of $A_U$; by the same token, $R$ contains
the $\bar A$-subalgebra $R'$ of $A_U$ generated by
$\phi^\flat_U(\bE_U)$, and proposition \ref{prop_boeuf}(ii)
implies that $R'$ is dense in $R$. However, $R'$ is
also open in $A_U$, so $R=R'$. Furthermore, a simple
inspection shows that the definition of $\phi^\flat_U$
depends only on the ring $A$ and the open subset $U$
(and is independent of the topology of $A$ or of $\bE$).
This shows that the image of $u_U$ is independent of
the topology of $\bE_U$.
Lastly, let $\cT_{\alpha_0}$ be the $\alpha_0$-adic
topology of $\bE$; it follows easily from example
\ref{ex_perfectoid}(i) and lemma \ref{lem_fontaine}
that the topological ring $(\bE,\cT_{\alpha_0})$ is
perfectoid, and $U_\bE$ contains the analytic locus
of $\Spec\,\bE$, relative to the topology $\cT_{\alpha_0}$,
whence the claim.
\end{pfclaim}

Henceforth, we assume that the topology of $\bE_U$
is $\alpha_0$-adic. Set $J_1:=\Ann_{\bar\bE}(\alpha_0)$,
let $J_2\subset\bar\bE$ be the radical of the ideal
$\alpha_0\bar\bE$, and $J_3:=J_1+J_2$; let
$\bar X_\bE:=\Spec\,\bar\bE$, $\bar X_A:=\Spec\,\bar A$,
and for $i=1,2,3$ set $\bar\bE_i:=\bar\bE/J_i$,
$\bar A_i:=W(\bar\bE_i)\otimes_{W(\bar A)}\bar A$,
$\bar X_{\bE,i}:=\Spec\,\bar\bE_i$ and
$\bar X_{A,i}:=\Spec\,\bar A_i$.
Lastly, let $j_{\bE,i}:\bar X_{\bE,i}\to\bar X_\bE$ and
$j_{A,i}:\bar X_{A,i}\to\bar X_A$ be the induced closed
immersions and set
$\bE_{U,i}:=j_{\bE,i*}\cO_{\bar X_{\bE,i}}(U_\bE)$,
$A_{U,i}:=j_{A,i*}\cO_{\bar X_{A,i}}(U_A)$ for $i=1,2,3$.
According to \eqref{subsec_decompose-torsion} we have
cartesian diagrams of quasi-coherent $\cO_{\bar X_\bE}$-modules
$$
\xymatrix{ \cO_{\bar X_\bE} \ar[r] \ar[d] &
j_{\,E,1*}\cO_{\bar X_{\bE,1}} \ar[d] &
\cO_{\bar X_A} \ar[r] \ar[d] &
j_{A,1*}\cO_{\bar X_{A,1}} \ar[d] \\
j_{\bE,2*}\cO_{\bar X_{\bE,2}} \ar[r] & j_{\bE,3*}\cO_{\bar X_{\bE,3}}
& j_{A,2*}\cO_{\bar X_{A,2}} \ar[r] & j_{A,3*}\cO_{\bar X_{A,3}}
}$$
whence cartesian diagrams of ring homomorphisms
\set\begin{equation}\label{eq_no-claim}
{\diagram \bE_U \ar[r] \ar[d] & \bE_{U,1} \ar[d]
& A_U \ar[r] \ar[d] & A_{U,1} \ar[d] \\
\bE_{U,2} \ar[r] & \bE_{U,3} & A_{U,2} \ar[r] & A_{U,3}.
\enddiagram}
\end{equation}
We endow each $\bE_{U,i}$ with its $\alpha_0$-adic topology
$\cT_i$; then $\cT_i$ is the discrete topology for
$i=2,3$, and the topological ring $(\bE_{U,i},\cT_i)$ is
perfectoid for $i=1,2,3$ (see
\eqref{subsec_decompose-torsion}). We also define the
power-multiplicative norms $|\cdot|_{I,i}:\bE_{U,i}\to\R_+$
associated with the ideal $I\bE_{U,i}$ as in
\eqref{eq_back-on-track}, and we get consequently the
normed topological rings $W(\bE_{U,i},1)$ as in
\eqref{subsec_boeuf}.

\begin{claim} The induced diagrams of rings
$$
\cW(\bE)\ :\ 
{\diagram W(\bE_U) \ar[r] \ar[d] & W(\bE_{U,1}) \ar[d] \\
W(\bE_{U,2}) \ar[r] & W(\bE_{U,3})
\enddiagram}
\qquad
\cW(\bE,1)\ :\ 
{\diagram W(\bE_U,1) \ar[r] \ar[d] & W(\bE_{U,1},1) \ar[d] \\
W(\bE_{U,2},1) \ar[r] & W(\bE_{U,3},1)
\enddiagram}
$$
are cartesian.
\end{claim}
\begin{pfclaim} The assertion is clear for the left
diagram. Next, set
$$
\bE^+_U:=\{x\in\bE_U~|~|x|_I\leq 1\}
$$
and define likewise the subring $\bE_{U,i}^+\subset\bE_{U,i}$
for $i=1,2$. For every $x\in\bE_U$ and $i=1,2$, denote by
$x^{(i)}\in\bE_{U,i}$ the image of $x$; since we already know
that diagram $\cW(\bE)$ is cartesian, the assertion for
$\cW(\bE,1)$ comes down to checking that if $x\in\bE_U$
and $|x^{(i)}|_{I,i}\leq\rho^{-s}$ for $i=1,2$ and some
$s\in\N[1/p]$, then $|x|_I\leq\rho^{-s}$. The latter is
equivalent in turn to the following. Suppose that for
every strictly positive $\eps\in\N[1/p]$ and every
$b\in I^{\La s+\eps\Ra}$ we have $b^{(i)}x^{(i)}\in\bar\bE_i$;
then $z:=bx\in\bE_U^+$. Now, for such $x$ and $b$ pick
$y_i\in\bar\bE$ with $y^{(i)}_i=b^{(i)}x^{(i)}$ for $i=1,2$;
then $(z-y_1)(z-y_2)\in J_1\cap J_2=0$, and especially,
$z$ is integral over $\bar\bE$. But $\bar\bE\subset\bE_U^+$,
and $\bE^+_U$ is integrally closed in $\bE_U$ (remark
\ref{rem_semi-norm}(vi)), whence the claim.
\end{pfclaim}

\begin{claim}\label{cl_tens-by-A}
The diagram of rings $\cW(\bE,1)\otimes_{W(\bE)}A$ is cartesian.
\end{claim}
\begin{pfclaim} Consider the complex in $\sC^{[0,1]}(W(\bE)\Mod)$
$$
K^\bullet
\quad :\quad
0\to W(\bE_{U,1},1)\oplus W(\bE_{U,2},1)\to W(\bE_{U,3},1)\to 0
$$
whose differential is deduced from $\cW(\bE,1)$, so that
$H^0K^\bullet=W(\bE_U,1)$. Taking into account proposition
\ref{prop_reduced-Witt}(i), we have a short exact sequence
$$
0\to K^\bullet\xrightarrow{\ \underline\alpha\cdot\one_{K^\bullet}\ }
K^\bullet\to K^\bullet/\underline\alpha K^\bullet\to 0
$$
and by the same token, scalar multiplication is injective
on $H^0K^\bullet$, so we get a long exact cohomology sequence
$$
0\to H^0(K^\bullet)/\underline\alpha\cdot H^0(K^\bullet)
\xrightarrow{\ f_{\underline\alpha}\ }
H^0(K^\bullet/\underline\alpha K^\bullet)\to
H^1K^\bullet\xrightarrow{\ \underline\alpha\ }H^1K^\bullet
$$
and the claim comes down to checking that $f_{\underline\alpha}$
is an isomorphism, or equivalently, that scalar multiplication
by $\underline\alpha$ is injective on $H^1K^\bullet$. However,
notice that $H^1K^\bullet$ is a quotient of the $W(\bE)$-module
$W(\bE_{U,3},1)$, and the latter is annihilated by
$\tau_\bE(\alpha_0)$.
Since $\underline\alpha=\tau_\bE(\alpha_0)+p\underline u$
for an invertible element $\underline u$ of $W(\bE)$, we
are reduced to checking that scalar multiplication by
$p$ is injective on $H^1K^\bullet$, {\em i.e.} we may assume
that $\underline\alpha=p$, and then it suffices to show
that $f_p$ is an isomorphism. However, it follows easily
from proposition \ref{prop_boeuf}(iii) that
\set\begin{equation}\label{eq_ALEAS}
pW(\bE_{U,i})\cap W(\bE_{U,i},1)=pW(\bE_{U,i},1)
\qquad
\text{for $i=1,2,3$}
\end{equation}
and likewise for $W(\bE_U,1)$. Therefore,
$H^0(K^\bullet)/p\cdot H^0(K^\bullet)=\bE_U$, and
$K^\bullet/pK^\bullet$ is the complex
$$
0\to\bE_{U,1}\oplus\bE_{U,2}\to\bE_{U,2}\to 0
$$
deduced from the cartesian diagram \eqref{eq_no-claim}.
The claim follows.
\end{pfclaim}

From claim \ref{cl_tens-by-A} we deduce a commutative diagram
of rings
$$
\xymatrix{ A_U \ar[rrr] \ar[ddd] & & & A_{U,i} \ar[ddd] \\
& W(\bE_U,1)\otimes_{W(\bE)}A \ar[r] \ar[d] \ar[lu]_{u_U} &
W(\bE_{U,1},1)\otimes_{W(\bE)}A \ar[d] \ar[ru]^{u_{U,1}} \\
& W(\bE_{U,2},1)\otimes_{W(\bE)}A \ar[r] \ar[ld]_{u_{U,2}}
& W(\bE_{U,3},1)\otimes_{W(\bE)}A \ar[rd]^{u_{U,3}} \\
A_{U,2} \ar[rrr] & & & A_{U,3}
}$$
whose two square subdiagrams are both cartesian. Recall
that the topologies of $\bE_{U,2}$ and $\bE_{U,3}$ are
discrete; taking into account theorem
\ref{th_int-subrings-perfectoid}(i) it follows that
$\bE_{U,i}=\bE^\circ_{U,i}=A^\circ_{U,i}=A_{U,i}$ and
$W(\bE_{U,i},1)=W(\bE_{U,i})$ for $i=2,3$, which easily
implies that $u_{U,2}$ and $u_{U,3}$ are isomorphisms.
Thus, to conclude the proof it suffices to show that
$u_{U,1}$ is an isomorphism. Hence, we may replace
$\bE$ by $\bar\bE_1$ and assume from start that
$\alpha_0$ is a regular element of $\bE$ and the
topology of $\bE$ is $\alpha_0$-adic. Let $R$ be
the image of $u_U$; under the current assumptions,
$p$ is a regular element of $A$ (remark
\ref{rem_decompose-torsion}(ii)), and therefore
$$
A\subset R\subset A_U\subset A[1/p]
$$
whence $R[1/p]=A[1/p]$. Moreover, we notice :

\begin{claim}\label{cl_final-claim}
The inclusion $R\subset A$ induces an injective map
$R/pR\to A_U/pA_U$.
\end{claim}
\begin{pfclaim} From (ii) and \eqref{eq_ALEAS} we get
a natural identification
\set\begin{equation}\label{eq_yes}
R/pR\isom W(\bE_U,1)/(\underline\alpha W(\bE_U,1)+pW(\bE_U,1))
\isom\bE_U/\alpha_0\bE_U
\end{equation}
and from the short exact sequence of $\cO_{\!X_\bE}$-modules
$$
0\to\cO_{\!X_\bE}\xrightarrow{\ \alpha_0\cdot\ }
\cO_{\!X_\bE}\to\cO_{\!X_\bE}/\alpha_0\cO_{\!X_\bE}\to 0
$$
we see that the natural map
$$
\bE_U/\alpha_0\bE_U\to\bE^0_U:=
H^0(U_\bE\cap\Spec\,\bE/\alpha_0\bE,\cO_{\!X_\bE}/\alpha_0\cO_{\!X_\bE})
$$
is injective. Likewise, we see that the same holds for
the natural map
$$
A_U/pA_U\to A^0_U:=H^0(U_A\cap\Spec\,A/pA,\cO_{\!X_A}/p\cO_{\!X_A}).
$$
On the other hand, recall that the ring homomorphism
$u_A:W(\bE)\to A$ induces an isomorphism of schemes
$\phi_0:\Spec\,A/pA\isom\Spec\,\bE/\alpha_0\bE$ (see
\eqref{subsec_phi-flat-map}), whence a ring isomorphism
$\phi^\flat_{0,U}:\bE^0_U\isom A^0_U$ and a simple inspection
of the definitions shows that the resulting diagram
$$
\xymatrix{
\bE_U/\alpha_0\bE_U \ar[rrr]^-{u_U\otimes_{W(\bE)}A/pA} \ar[d]
& & & A_U/pA_U \ar[d] \\
\bE^0_U \ar[rrr]^-{\phi^\flat_{0,U}} & & & A^0_U
}$$
commutes. But it is easily seen that the isomorphism
\eqref{eq_yes} identifies $u_U\otimes_{W(\bE)}A/pA$ with
the map $R/pR\to A_U/pA_U$ induced by the inclusion
$R\subset A$, whence the claim.
\end{pfclaim}

Now, let $a\in A_U$ be any element, and pick $n\in\N$
such that $p^na\in A$; then there exists
$\underline x\in W(\bE)$ such that
$u_U(\underline x)=u_A(\underline x)=p^na$. To conclude,
we shall show, by induction on $i=0,\dots,n$ that
$p^{n-i}a\in R$. Indeed, we have just shown that $p^na\in R$.
Suppose that $i\geq 0$ and $p^{n-i}a\in R$; if $i=n$,
we are done, and if $i<n$, the image of $p^{n-i}a$
vanishes in $A_U/pA_U$, hence also in $R/pR$, {\em i.e.}
there exists $b\in R$ with $pb=p^{n-i}a$; since $p$ is
regular in $R$, it follows that $b=p^{n-i-1}a$, as required.
\end{proof}

\begin{corollary}\label{cor_boody-boody}
In the situation of \eqref{subsec_boeuf}, we have :
\begin{enumerate}
\item
$\{\Fil^s\bE_U\}=\Fil^sA_U$ for every $s\in\Z[1/p]$.
\item
The map $\phi^\flat_U$ induces a ring isomorphism
$\bE_U/\alpha_0\bE_U\isom A_U/pA_U$.
\end{enumerate}
\end{corollary}
\begin{proof}(i): Let us first remark that
$$
\{\Fil^{s-n}\bE_U\}\cap\Fil^sA_U=\{\Fil^s\bE_U\}
\qquad
\text{for every $n\in\N$ and every $s\in\Z[1/p]$}.
$$
Indeed, the case where $n=1$ follows immediately from
claim \ref{cl_spec-case-qaff}(ii), and then the general
case follows by a simple induction on $n\in\N$ : details
left to the reader. But theorem \ref{th_boody-boody}(i)
shows that $A_U=\bigcup_{n\in\N}\{\Fil^{s-n}\bE_U\}$, whence
the assertion.

(ii) follows directly from theorem \ref{th_boody-boody}(i)
and \eqref{eq_yes}.
\end{proof}

As an application, we point out the following criterion :

\begin{corollary}\label{cor_E-preserves-affines}
Let $A$ be any perfectoid ring,
$U_\bE\subset V_\bE\subset X_\bE:=\Spec\,\bE(A)$ two
constructible open subsets that contain the analytic
locus. Set also $U_A:=\phi^{-1}U_\bE$ and $V_A:=\phi^{-1}V_\bE$,
where $\phi:X_A:=\Spec\,A\to X_\bE$ is the continuous
map of \eqref{subsec_map-on-spec}. Then we have :
\begin{enumerate}
\item
$U_\bE$ is an affine scheme if and only if the same
holds for $U_A$.
\item
More generally, the inclusion map $U_\bE\to V_\bE$ is
an affine morphism of schemes if and only if the same
holds for the inclusion map $U_A\to V_A$.
\end{enumerate}
\end{corollary}
\begin{proof}(i): Set $\bE_U:=\cO_{\!X_\bE}(U_\bE)$ and
$A_U:=\cO_{\!X_A}(U_A)$; by lemma \ref{lem_dots-on-is}(i),
we get commutative diagrams of schemes
$$
\xymatrix{ U_\bE \ar[r]^-{j_\bE} \ar[rd]_{i_\bE} &
\Spec\,\bE_U \ar[d]^{j'_\bE} &
U_A \ar[r]^-{j_A} \ar[rd]_{i_A} & \Spec\,A_U \ar[d]^{j'_A}
\\
& X_\bE & & X_A
}$$
where $i_\bE$, $i_A$, $j_\bE$ and $j_A$ are open immersions.
We have to show that $j_\bE$ is an isomorphism if and only
if the same holds for $j_A$. To this aim, set
$Z_\bE:=X_\bE\setminus U_\bE$ and $Z_A:=X_A\setminus U_A$;
equivalently, we have to check that $j'^{-1}_\bE Z_\bE=\emptyset$
if and only if $j'^{-1}_AZ_A=\emptyset$. Now, since $U_\bE$
contains the analytic locus of $X_\bE$, we may assume that
$Z_\bE\subset X_\bE^0:=\Spec\,\bE/\alpha_0\bE$, where
$(\alpha_n~|~n\in\N)$ is any fixed distinguished element
of $\Ker\,u_A$, and likewise $Z_A\subset X_A^0:=\Spec\,A/pA$.
On the other hand, by corollary \ref{cor_boody-boody}(ii)
we have a commutative diagram of schemes
$$
\xymatrix{ \Spec\,A_U/pA_U \ar[r] \ar[d] &
\Spec\,\bE_U/\alpha_O\bE_U \ar[d] \\
X_A^0 \ar[r]^-{\phi_0} & X_\bE^0
}$$
whose horizontal arrows are isomorphisms induced by $\phi$
and by the map $\phi_U^\flat:\bE_U\to A_U$ of proposition
\ref{prop_new-formula}(i), and whose vertical arrows are the
restrictions of $j'_\bE$ and $j'_A$. Lastly, by construction
we have $\phi_0^{-1}Z_\bE=Z_A$, whence the contention.

(ii) is analogous : the open immersion
$U_\bE\to V_\bE$ is affine if and only if for every affine
open subset $W\subset V_\bE$, the intersection $W\cap U_\bE$
is still affine. The latter condition means that
$$
\Spec\,\cO_{\!X_\bE}(W\cap U_\bE)\times_WZ_\bE=\emptyset
\qquad
\text{for every such $W$}
$$
which in turns is equivalent to the condition :
\set\begin{equation}\label{eq_recondition}
j_\bE'^{-1}(Z_\bE\cap V_\bE)=\emptyset.
\end{equation}
Likewise, we easily see that the inclusion $U_A\subset V_A$
is affine if and only if $j_A'^{-1}(Z_A\cap V_A)=\emptyset$,
and since $Z_\bE\subset X^0_\bE$, this latter condition is
equivalent to \eqref{eq_recondition} : details left to the reader.
\end{proof}

\sset\subsubsection{}\label{subsec_powers-and-Witt}
Let $\cK,\cJ\subset\bE_U$ be any two $\bar\bE$-submodules;
we let $\cK\cJ\subset\bE_U$ be the $\bar\bE$-submodule
generated by the system $(xy~|~x\in\cK,y\in\cJ)$. Especially,
we have well defined $\bar\bE$-submodules $\cK^n$ of $\bE_U$,
defined inductively by the rule : $\cK^0:=\bar\bE$ and
$\cK^{n+1}:=\cK^n\cK$ for every $n\in\N$. With this notation,
we also set
$$
W(\cK):=
\{(x_n~|~n\in\N)~|~x_n\in\cK^{p^n}\ \text{for every $n\in\N$}\}
\qquad
W(\cK,1):=W(\cK)\cap W(\bE_U,1).
$$
Just as in remark \ref{rem_Witt-limit}, it is easily seen
that $W(\cK)$ and $W(\cK,1)$ are $W(\bar\bE)$-submodule of
respectively $W(\bE_U)$ and $W(\bE_U,1)$. Denote also by
$(\cK^n)^c$ the topological closure of $\cK^n$ in $\bE_U$,
for every $n\in\N$, and by $W(\cK,1)^c$ the topological
closure of $W(\cK,1)$ in $W(\bE_U,1)$. 

\begin{lemma}\label{lem_like-th-taut-two}
With the notation of \eqref{subsec_powers-and-Witt}, we have
\begin{enumerate}
\item
$W(\cK,1)^c=M:=\{(x_n~|~n\in\N)\in W(\bE_U,1)~|~
x_n\in(\cK^{p^n})^c\ \text{for every $n\in\N$}\}$.
\item
If $\cK$ is taut, $u_U(W(\cK,1)^c)=\{\cK\}$.
\item
If $\cK$ is $1$-taut and topologically closed in $\bE_U$, then
$u_U$ induces an isomorphism
$$
W(\cK,1)/\underline\alpha W(\cK,1)\isom\{\cK\}.
$$
\end{enumerate}
\end{lemma}
\begin{proof}(i): Let $\underline x:=(x_n~|~n\in\N)\in M$ be any
element, and fix any $\eps\in\R_{>0}$. We shall exhibit inductively
a sequence $\underline y:=(y_n~|~n\in\N)$ of elements of $\bE_U$
such that
$$
|\underline y|_1\leq\eps
\qquad\qquad\text{and}\qquad
\underline x+\underline y\in W(\cK,1).
$$
Indeed, let $n\in\N$, and suppose that we have already
exhibited $y_i\in(\cK^{p^i})^c$ for $i=0,\dots,n-1$ such that
$\underline y^{(n)}:=\sum_{i=0}^{n-1}p^i\cdot\tau_{\bE_U}(y_i)$ has
norm $<\eps$ and such that if we set
$(z^{(n)}_i~|~i\in\N):=\underline x+\underline y^{(n)}$, we have
$z^{(n)}_i\in\cK^{p^i}$ for every $i=0,\dots,n-1$. It follows
easily that $z^{(n)}_i\in(\cK^{p^i})^c$ for every $i\in\N$. We
may then pick $y_n\in(\cK^{p^n})^c$ such that $|y_n|_I<\eps^{p^n}$
and $z^{(n)}_n+y_n\in\cK^{p^n}$. Set
$\underline y^{(n+1)}:=\underline y^{(n)}+p^n\cdot\tau_{\bE_U}(y_n)$
and $(z^{(n+1)}_i~|~i\in\N):=\underline x+\underline y^{(n+1)}$.
In light of lemma \ref{lem_Witt-truncate}(i) and claim
\ref{cl_W_n-structure} it is easily seen that
$z^{(n+1)}_i=z^{(n)}_i$ for every $i=0,\dots,n-1$, and
$z^{(n+1)}_n=z^{(n)}_n+y_n$. Then, it is clear that the element
$\underline y:=\sum_{n\in\N}p^n\cdot\tau_{\bE_U}(y_n)$ will do.

This shows that $M\subset W(\cK,1)^c$. The converse inclusion
is obvious, since $\prod_{n\in\N}(\cK^{p^n})^c$ is a closed subset
of $W(\bE_U)$ that contains $W(\cK,1)$ and the topology of
$W(\bE_U,1)$ is finer than the one induced by $W(\bE_U)$.

(ii): Recall that
$\phi^\flat_U(\alpha_0)=\bar u_{\bar A}(\alpha_0)=pv$ for an
element $v\in\bar A{}^\times$. Then, since $\cK$ is taut,
we see that $p^n\cdot\phi_U^\flat(x^{1/p^n})=
v^{-n}\phi^\flat_U(\alpha_0^nx^{1/p^n})\in\{\cK\}$ for every
$x\in\cK^{p^n}$. It follows easily that
$u_U(\underline x)\in\{\cK\}$ for every
$\underline x\in W(\cK)$, and thus
$u_U(W(\cK,1)^c)\subset\{\cK\}$.

For the converse inclusion, pick a strictly positive
$\eps\in\N[1/p]$ such that $\alpha_0\in I^{\La\eps\Ra}$; it was
observed in the proof of lemma \ref{lem_boody-boody}(iii)
that $\phi_U^\flat$ induces an $\bar A/p\bar A$-linear
isomorphism
$$
(\Fil^s\cK)/(\Fil^{s+\eps}\cK)\isom
(\Fil^s\{\cK\})/(\Fil^{s+\eps}\{\cK\})
\qquad
\text{for every $s\in\Z[1/p]$}
$$
and moreover $\Fil^s\{\cK\}=\{\cK\}\cap\{\Fil^s\bE_U\}=
\{\cK\}\cap\Fil^sA_U$ for every $s\in A_U$ (lemma
\ref{lem_boody-boody}(iii) and corollary \ref{cor_boody-boody}),
hence the filtration $\Fil^\bullet\{\cK\}$ is exhaustive and
separated on $\{\cK\}$. Thus, say that $x\in\Fil^s\{\cK\}$
for some $s\in\Z[1/p]$; it follows easily that we may find
a sequence $(y_n~|~n\in\N)$ of elements of $\bE_U$ such
that $y_n\in\Fil^{s+n\eps}\cK$ for every $n\in\N$, and the
series $\sum_{n\in\N}\phi_U^\flat(y_n)$ converges to $x$ in
the topology of $A_U$. Clearly $\tau_{\bE_U}(y_n)\in W(\cK)$
and $|\tau_{\bE_U}(y_n)|_1\leq\rho^{s+n\eps}$ for every $n\in\N$,
so the series $\sum_{n\in\N}\tau_{\bE_U}(y_n)$ converges in
$W(\bE_U,1)$ to an element $\underline z\in W(\cK)^c$,
and finally
$u_U(\underline z)=\sum_{n\in\N}u_U\circ\tau_{\bE_U}(y_n)=x$.

(iii): The assumptions on $\cK$ imply that $\cK^{p^n}$ is
topologically closed in $\bE_U$, for every $n\in\N$;
combining with (i) and (ii), it follows already that the
induced map $W(\cK,1)/\underline\alpha W(\cK,1)\to\{\cK\}$
is surjective. Next, say that $\alpha_0\in I^{\La\eps\Ra}$ for
some element $\eps>0$ of $\Z[1/p]$; notice that $\Fil^s\cK$
is still $1$-taut and topologically closed in $\bE_U$ for
every $s\in\Z[1/p]$; taking into account lemma
\ref{lem_boody-boody}(iii), we may therefore reduce to
the case where $\cK\subset\Fil^{-n\eps}\bE_U$ for some $n\in\N$.
We argue by induction on $n$ : indeed, the assertion for
$n=0$ is already known, by theorem \ref{th_taut-two}(iv).
Thus, suppose that $n>0$, and that the assertion is already
known for $\cK\cap\Fil^{(1-n)\eps}\bE_U$; let
$\underline x:=(x_n~|~n\in\N)\in W(\cK,1)$ be an element
in the kernel of $u_U$. This means that
$\sum_{n\in\N}p^i\cdot\phi^\flat_U(x^{1/p^n}_n)=0$ in $A_U$, and we
need to check that $\underline x\in\underline\alpha W(\cK,1)$.
Now, since $\cK$ is $1$-taut, we have
$\underline z:=(x_{n+1}^{1/p}~|~n\in\N)\in W(\cK,1)$ as well,
so that
$$
\phi^\flat_U(x_0)=
-\sum_{n\in\N}p^{n+1}\cdot\phi^\flat_U(x_{n+1}^{1/p^{n+1}})
\in\Fil^{(1-n)\eps}A_U.
$$
In view of claim \ref{cl_spec-case-qaff}(ii), we deduce
that $x_0\in\Fil^{(1-n)\eps}\cK$. Moreover,
$\underline x=\tau_{\bE_U}(x_0)+p\cdot\underline z$, and recall
that $\underline\alpha=\tau_\bE(\alpha_0)+p\cdot\underline u$
for an invertible element $\underline u\in W(\bE)$; we
conclude that $\underline x$ is the sum of
$$
\tau_{\bE_U}(x_0)-
\underline u^{-1}\cdot\tau_\bE(\alpha_0)\cdot\underline z
\in W(\Fil^{(1-n)\eps}\cK,1)
\qquad\text{and}\qquad
\underline u^{-1}\cdot\underline\alpha\cdot\underline z
\in\underline\alpha W(\cK,1).
$$
Especially, $\tau_{\bE_U}(x_0)-
\underline u^{-1}\cdot\tau_\bE(\alpha_0)\cdot\underline z$
lies as well in the kernel of $u_U$; by inductive assumption
it is therefore an element of
$\underline\alpha W(\Fil^{(1-n)\eps}\cK,1)$, and the proof is
concluded.
\end{proof}

\begin{proposition}\label{prop_Cont-u_A}
Let $A$ be any complete and separated topological ring,
whose topology is linear and coarser than the $p$-adic
topology, and $v\in\Cont(A)$ any valuation. We have :
\begin{enumerate}
\item
The mapping $v\circ\bar u_A$ is a continuous valuation of
the topological ring $\bE:=\bE(A)$ (see remark
{\em\ref{rem_fontaine}(ii)}).
\item
The rule : $v\mapsto v\circ\bar u_A$ defines continuous
maps
$$
\Cont(\bar u_A):\Cont(A)\to\Cont(\bE)
\qquad
\Cont^+(\bar u_A):\Cont^+(A)\to\Cont^+(\bE).
$$
\item
The diagrams of continuous maps
$$
\xymatrix{ \Cont(A) \ar[rr]^-{\Cont(\bar u_A)} \ar[d] & &
\Cont(\bE) \ar[d] & & \Cont^+(A) \ar[rr]^-{\Cont^+(\bar u_A)}
\ar[d] & & \Cont^+(\bE) \ar[d] \\
\Spec\,A \ar[rr]^-{\Spec\,\bar u_A} & & \Spec\,\bE
& & \Spec\,A \ar[rr]^-{\Spec\,\bar u_A} & & \Spec\,\bE
}$$
commute, where $\Spec\,\bar u_A$ is defined as in
\eqref{sec_def-A-tilde} and the vertical arrows
in the left (resp. right) diagram are the restrictions
of the support maps $\sigma_A$ and $\sigma_\bE$ (resp.
of the center maps $\sigma_A^+$ and $\sigma^+_\bE$) : see
remark {\em\ref{rem_Spv-of-ring}(iii,v)}.
\item
If\/ $v'\in\Cont(A)$ is a primary (resp. secondary) specialization
of\/ $v$, then $\Cont(\bar u_A)(v')$ is a primary (resp. secondary)
specialization of\/ $\Cont(\bar u_A)(v)$.
\item
The map $\Cont(\bar u_A)$ restricts to a surjection from
the set of secondary specializations of\/ $v$ onto the set
of secondary specializations of\/ $\Cont(\bar u_A)(v)$.
\item
If $A$ is a P-ring, $\Cont(\bar u_A)$ also restricts to a
surjection from the set of secondary generizations of\/ $v$
onto the set of secondary generizations of\/ $\Cont(\bar u_A)(v)$.
\end{enumerate}
\end{proposition}
\begin{proof}(i): Since both $v$ and $\bar u_A$ are continuous
morphisms of pointed monoids, the same holds for $v\circ\bar u_A$.
Thus, it remains only to check that
$$
v(\bar u_A(a_\bullet+b_\bullet))\leq
\max(v(a_0),v(b_0))
\qquad
\text{for every $a_\bullet,b_\bullet\in\bE$}.
$$
However, recall that
$$
\bar u_A(a_\bullet+b_\bullet)=\lim_{n\to+\infty}(a_n+b_n)^{p^n}
$$
(remark \ref{rem_fontaine}(ii)). Since the pointed
monoid $\Gamma_{\!v\circ}$ is separated for its topology
$\cT_{\Gamma_{\!v}}$ (notation of definition
\ref{def_continuous-vals}(i)), it follows that
$v(\bar u_A(a_\bullet+b_\bullet))$ is the unique limit
point of the sequence $(v(a_n+b_n)^{p^n}~|~n\in\N)$.
However :
$$
v(a_n+b_n)^{p^n}\leq\max(v(a_n)^{p^n},v(b_n)^{p^n})=
\max(v(a_0),v(b_0))
$$
whence the assertion.

(ii): Let $a_\bullet,b_\bullet\in\bE$ be any two elements,
and set
$U:=\{v\in\Cont(\bE)~|~v(a_\bullet)\leq v(b_\bullet)\neq 0\}$;
directly from the definition we find :
$$
\Cont(\bar u_A)^{-1}U=\{w\in\Cont(A)~|~w(a_0)\leq w(b_0)\neq 0\}
$$
whence the assertion.

(iii) and (iv) are immediate from the definitions.

(v): Set $w:=\Cont(\bar u_A)(v)$ and notice that $\fp:=\sigma_A(v)$
is a closed subset in the $p$-adic topology of $A$; we remark
the following :

\begin{claim} Denote by $\bar v$ (resp. $\bar w$) the
residual valuation of $v$ (resp. of $w$). Then we have :
\begin{enumerate}
\item
The map $\bar u_\fp:\kappa(w)\to\kappa(v)$ of lemma
\ref{lem_new-place} restricts to a local morphism
of multiplicative monoids
$\bar u{}^+_v:(\kappa(w),\bar w)^+\to(\kappa(v),\bar v)^+$.
\item
Let also $\bar\kappa(v)$ (resp. $\bar\kappa(w)$) be the
residue field of $\kappa(v)^+$ (resp. of $\kappa(w)^+$).
There exists a unique ring homomorphism
$\bar\phi_v:\bar\kappa(w)\to\bar\kappa(v)$ fitting
into a commutative diagram
$$
\xymatrix{
\kappa(w)^+ \ar[rr]^-{\bar u^+_v} \ar[d] & &
\kappa(v)^+ \ar[d] \\
\bar\kappa(w) \ar[rr]^-{\bar\phi_v} & & \bar\kappa(v)
}$$
whose vertical arrows are the projections.
\end{enumerate}
\end{claim}
\begin{pfclaim}(i): Indeed, a simple inspection of the
definitions yields a commutative diagram of monoids
$$
\xymatrix{
\kappa(w) \ar[rr]^-{\bar u_\fp} \ar[rd]_{\bar w} & &
\kappa(v) \ar[ld]^{\bar v} \\
& \Gamma_{\!v}
}$$
whence the assertion : details left to the reader.

(ii): Let $\bar\pi_v:\kappa(v)^+\to\bar\kappa(v)$ be the
projection; since $\bar u^+_v$ is local, we are easily
reduced to checking that the composition
$\psi:=\bar\pi_v\circ\bar u^+_v$ is a ring homomorphism.
However, $\psi$ is obviously a morphism of multiplicative
monoids, so it suffices to show that
$$
\bar v(\bar u^+_v(x_1+x_2)-\bar u^+_v(x_1)-\bar u^+_v(x_2))<1
\qquad
\text{for every $x_1,x_2\in\kappa(w)^+$}.
$$
However, let $\pi_w:\bE\to\kappa(w)$ and $\pi_v:A\to\kappa(v)$
be the projections; we may write
$$
x_1=\pi_w(\beta_1)/\pi_w(\gamma)
\qquad\text{and}\qquad
x_2=\pi_w(\beta_2)/\pi_w(\gamma)
\qquad
\text{for some $\beta_1,\beta_2,\gamma\in\bE$}
$$
with $w(\beta_1),w(\beta_2)\leq w(\gamma)\neq 0$.
Set $x_3:=x_1+x_2$ and $\beta_3:=\beta_1+\beta_2$; then
$$
\bar u^+_v(x_i)=
\pi_v\circ\bar u_A(\beta_i)/\pi_v\circ\bar u_A(\gamma)
\qquad
\text{for $i=1,2,3$}.
$$
Let $y:=\bar u_A(\beta_3)-\bar u_A(\beta_1)-\bar u_A(\beta_2)$;
we are then reduced to checking that $v(y)<v(\bar u_A(\gamma))$.
However, proposition \ref{prop_combinatorial} says that
$y=\sum_{n\in\N}p^{n+1}z_n$ where each $z_n$ is a finite
$\Z_p$-linear combination of terms of the form
$\bar u_A(\beta_1^\lambda\beta_2^{1-\lambda})$ with
$\lambda,1-\lambda\in\N[1/p]$. Clearly
$w(\beta_1^\lambda\beta_2^{1-\lambda})\leq w(\gamma)$ for
every such $\lambda$, so that
$$
v(p^{n+1}z_n)\leq v(p)^{n+1}\cdot v(\bar u_A(\gamma))
\qquad
\text{for every $n\in\N$}.
$$
Since $v$ is a continuous valuation and $v(p)$ is
final in $\Gamma_{\!v}$, the contention follows
easily : details left to the reader.
\end{pfclaim}

Now we argue as in the proof of lemma
\ref{lem_image-of-special}(iii). Namely, let $w'$ be a
secondary specialization of $w$; then $\kappa(w')=\kappa(w)$,
and $\kappa(w')^+\subset\kappa(w)^+$. The image $W$ of
$\kappa(w')^+$ in $\bar\kappa(w)$ is a valuation ring, and
by corollary \ref{cor_cornerstone} there exists a valuation
ring $V$ of $\bar\kappa(v)$ that dominates $\bar\phi_v(W)$.
The preimage of $V$ in $\kappa(v)^+$ is a valuation ring
of $\kappa(v)$, and the corresponding valuation $v'$ is
a secondary specialization of $v$. Lastly, $v'$ is continuous,
by remark \ref{rem_Cont-A}(iv), and by construction we have
$\Cont(\bar u_A)(v')=w'$.

(vi): Clearly $\bar u_A$ induces an injective morphism of
ordered groups $i:\Gamma_{\!w}\to\Gamma_{\!v}$, and $\Spec\,i$
is surjective, by remark \ref{rem_ordered-gps}(v). Now, let
$w'$ be a secondary generization of $w:=\Cont(\bar u_A)(v)$,
so that $w'=w_\Delta$ for a convex subgroup
$\Delta\subset\Gamma_{\!w}$, and we may write
$\Delta=i^{-1}\Delta'$ for some convex subgroup
$\Delta'\subset\Gamma_{\!v}$. If $\Delta'\neq\Gamma_{\!v'}$,
the projection $\Gamma_{\!v\circ}\to(\Gamma_{\!v}/\Delta')_\circ$
is continuous (remark \ref{rem_Cont-A}(iii)), so $v':=v_{\Delta'}$
is a secondary specialization of $v$ in $\Cont(A)$ with
$\Cont(\bar u_A)(v')=w'$. Hence, we may assume that $w'$
is a trivial valuation. It follows easily that
$\bE^{\circ\circ}\subset\Ker\,w'$. On the other hand, we may
find finitely many elements
$\beta_1,\dots,\beta_k\in\bE^{\circ\circ}$ such that the system
$\bar u_A(\beta_1),\dots,\bar u_A(\beta_k)$ generates an open
ideal of $A$ (see the proof of lemma \ref{lem_was-third-cond}(i)).
Summing up, we see that $v'$ is an analytic valuation of $A$,
especially $v'\in\Cont(A)$, and the proof is concluded.
\end{proof}

\begin{remark}\label{rem_Cont-u_A}
Keep the notation of proposition \ref{prop_Cont-u_A};
for every $v\in\Cont(A)$ we have the commutative diagram
of topological monoids
$$
\xymatrix{
\bE \ar[r]^-{\bar u_A} \ar[d]_{\bE(v)} & A \ar[d]^v \\
\bE(\Gamma_{\!v\circ}) \ar[r]^-{\bar u_{\Gamma_{\!v\circ}}} &
\Gamma_{\!v\circ}
}$$
(where $\Gamma_{\!v\circ}$ is endowed with its natural topology
$\cT_{\Gamma_{\!v}}$ as in definition \ref{def_continuous-vals}(i))
and under the natural identification
$\bE(\Gamma_{\!v\circ})\isom\bE(\Gamma_{\!v})_\circ$,
the map $\bar u_{\Gamma\!v\circ}$ corresponds to
$(\bar u_{\Gamma_{\!v}})_\circ$. Then $\bE(\Gamma_{\!v})$
inherits from $\Gamma_{\!v}$ a unique ordering such
that $\bar u_{\Gamma_{\!v}}$ is an injective morphism of
ordered abelian groups. Moreover, since the $p$-Frobenius
endomorphism of $\Gamma_{\!v}$ is injective, it is easily seen
that the topology of $\bE(\Gamma_{\!v\circ},\cT_{\Gamma_{\!v}})$
is the same as the topology $\cT_{\bE(\Gamma_{\!v})}$ provided
by definition \ref{def_continuous-vals}(i)  (details left
to the reader). In view of proposition \ref{prop_Cont-u_A}(i),
it follows that $\bE(v)$ is a continuous valuation on $\bE(v)$,
and we have
$$
\Cont(\bar u_A)(v)=\bE(v)
\qquad
\text{in $\Cont(\bE)$}.
$$
\end{remark}

\begin{proposition}\label{prop_A-is-valring-is-E}
Let $A$ be any perfectoid ring. We have :
\begin{enumerate}
\item
$A$ is a valuation ring if and only if the same holds
for $\bE:=\bE(A)$.
\item
Suppose that $A$ is a valuation ring, and let
$$
v_A:A\to\Gamma_{\!A\circ}:=(\Frac\,A)/A^\times
\qquad\text{and}\qquad
v_\bE:\bE\to\Gamma_{\!\bE\circ}:=(\Frac\,\bE)/\bE^\times
$$
be the valuations of $A$ and $\bE$ (see remark
{\em\ref{rem_valuations}(iv)}). Then :
\begin{enumerate}
\item
$v_A\in\Cont(A)$ and $v_\bE\in\Cont(\bE)$.
\item
The map $\bar u_A$ induces an isomorphism of ordered
abelian groups :
$$
\gamma_A:\Gamma_{\!\bE}\isom\Gamma_{\!A}
$$
fitting into a commutative diagram
$$
\xymatrix{
\bE \ar[rr]^-{\bar u_A} \ar[d]_{v_\bE} & &
A \ar[d]^{v_A} \\
\Gamma_{\!\bE\circ} \ar[rr]^-{\gamma_{A\circ}} & & \Gamma_{\!A\circ}.
}$$
\item
More precisely, if $(\alpha_n~|~n\in\N)$ is any
distinguished element of $\Ker\,u_A$, the map
$\bar u_A$ induces a group isomorphism
$$
\Frac(\bE)^\times/(1+\alpha_0\bE)\isom\Frac(A)^\times/(1+pA).
$$
\end{enumerate}
\item
Let\/ $V$ be a valuation ring, and $\cT_V$ a linear,
complete, separated and f-adic topology on $V$ that
is coarser than the $p$-adic topology. Suppose that :
\begin{enumerate}
\item
The value group of\/ $V$ is $p$-divisible.
\item
The Frobenius endomorphism of\/ $V/pV$ is surjective.
\end{enumerate}
Then $(V,\cT_V)$ is a perfectoid ring.
\end{enumerate}
\end{proposition}
\begin{proof}(i): We remark :

\begin{claim}\label{cl_if-E-is-valring}
If $\bE$ is a valuation ring, then for every $a\in A$
there exist $t\in A^\times$ and $e\in\bE$ such that
$a=t\cdot\bar u_A(e)$.
\end{claim}
\begin{pfclaim} Fix any distinguished element
$(\alpha_n~|~n\in\N)$ of $\Ker\,u_A$, and recall that
$p=w\cdot\bar u_A(\alpha_0)$ for some $w\in A^\times$
(lemma \ref{lem_was-third-cond}(iii)). Since the $p$-adic
topology is separated on $A$, there exists $n\in\N$ such
that $a\in p^nA\setminus p^{n+1}A$, so that $a=p^nb$ for some
$b\in A\setminus pA$. Then there exists $\beta\in\bE$
such that $b-\bar u_A(\beta)=pz$ for some $z\in A$
(lemma \ref{lem_was-third-cond}(i)), and notice that
$\beta\notin\alpha_0\bE$, since otherwise we would have
$\bar u_A(\beta)\in pA$, whence $b\in pA$, a contradiction.
Thus, $\alpha_0=\beta\cdot\gamma$ for some element $\gamma$
of the maximal ideal of $\bE$, and we may write
$$
x=p^nb=p^n\cdot(\bar u_A(\beta)+pz)=w^n\cdot
\bar u_A(\alpha_0^n\beta)\cdot(1+\bar u_A(\gamma)\cdot wz).
$$
We need to check that $u:=1+\bar u_A(\gamma)\cdot wz$
is invertible in $A$; to this aim, since $p$ lies in the
Jacobson radical of $A$, it suffices to show that
the image of $u$ is invertible in $A/pA$, and we
are further reduced to checking that the image $\bar\gamma$
of $\bar u_A(\gamma)$ lies in the Jacobson radical
of $A/pA$. However, $u_A$ induces a ring isomorphism
$\omega:\bE/\alpha_0\bE\isom A/pA$ (remark
\ref{rem_nice-topology}(ii)), and $\omega^{-1}(\bar\gamma)$
is the class of $\gamma$, whence the contention.
\end{pfclaim}

Now, suppose that $\bE$ is a valuation ring;
directly from claim \ref{cl_if-E-is-valring} and
corollary \ref{cor_perf-are-reduced}(ii) we deduce that
$A$ is a domain. Then, let $a_1,a_2\in A$ be any two
elements, and write $a_i=t_i\cdot\bar u_A(e_i)$
with $t_i\in A^\times$ and $e_i\in\bE$ for $i=1,2$;
we may assume that $e_1\in e_2\bE$, whence $a_1\in a_2A$,
which easily implies that $A$ is a valuation ring.

Conversely, suppose that $A$ is a valuation ring; if the
field of fractions $K$ of $A$ has characteristic $0$,
we know already from lemma \ref{lem_E-is-val.ring}
that $\bE$ is a valuation ring. If the characteristic
of $K$ equals $p$, the ring $\bE$ is isomorphic to $A$,
whence the assertion.

(ii.a): In light of proposition \ref{prop_stays-valuation}(i,ii)
we see that either the topology $\cT_A$ of $A$ is discrete, or else
it agrees with the valuation topology of $A$. In either case, it
is clear that $v_A\in\Cont(A)$. The same applies to $\bE$,
since the latter is perfectoid as well.

(ii.b): From claim \ref{cl_if-E-is-valring} we immediately
see that $\gamma_A$ is surjective, and the injectivity follows
from remark \ref{rem_Cont-u_A}.

(ii.c): By (ii.b) and the snake lemma, we are reduced
to showing that $\bar u_A$ induces a group isomorphism
$$
\bE^\times/(1+\alpha_0\bE)\isom A^\times/(1+pA).
$$
However, $\bar u_A$ induces a ring isomorphism
$\bE/\alpha_0\bE\isom A/pA$ (remark
\ref{rem_nice-topology}(ii)), whence a group isomorphism
$(\bE/\alpha_0\bE)^\times\isom(A/pA)^\times$, and it remains
only to observe that the natural maps
$$
\bE^\times/(1+\alpha_0\bE)\to(\bE/\alpha_0\bE)^\times
\qquad
A^\times/(1+pA)\to(A/pA)^\times
$$
are isomorphisms, since $\alpha_0\bE$ (resp. $pA$) lies
in the Jacobson radical of $\bE$ (resp. of $A$).

(iii): In light of proposition \ref{prop_stays-valuation}(i),
we see that $\cT_V$ is either the discrete topology or the
valuation topology on $V$. If $\cT_V$ is discrete and
coarser than the $p$-adic topology, then $V$ must be
an $\F_p$-algebra, so in this case the assertion is a
special case of example \ref{ex_perfectoid}(i).
If $\cT_V$ is the valuation topology, let us endow
$K:=\Frac\,V$ with its valuation topology $\cT_K$; then
$(K,\cT_K)$ is a Tate ring, by proposition
\ref{prop_stays-valuation}(ii), hence $V$ admits an
ideal of adic definition of the form $I:=aV$, for a
non-zero topologically nilpotent element $a\in V$
(corollary \ref{cor_Tate}(ii)). Since the value group
of $V$ is $p$-divisible, and $p$ is topologically
nilpotent in $V$, we may assume that $p\in I^p$, in
which case, assumption (iii.b) implies that the
Frobenius endomorphism of $V/pV$ induces an isomorphism
$V/I\isom V/I^p$. By the same token, we easily see that
$(V,\cT_V)$ is a P-ring, so the assertion follows from
theorem \ref{th_regular-seq-criterion}.
\end{proof}

\begin{remark}\label{rem_A-is-valring-is-E}
Let $A$, $A'$ be two perfectoid valuation rings with
valuations $v:A\to\Gamma_{\!A\circ}$,
$v':A'\to\Gamma_{\!A'\circ}$, and $f:A\to A'$ a continuous ring
homomorphism; according to proposition
\ref{prop_A-is-valring-is-E}(i), both $\bE:=\bE(A)$ and
$\bE':=\bE(A')$ are valuation rings,
and we denote by $v_\bE:\bE\to\Gamma_{\!\bE\circ}$,
$v_{\bE'}:\bE'\to\Gamma_{\!\bE'\circ}$ their respective valuations.
Then it follows easily from proposition
\ref{prop_A-is-valring-is-E}(ii.c) that $f$ is injective if
and only if the same holds for $\bE(f):\bE\to\bE'$. Moreover,
if $f$ is injective there exist unique group homomorphisms
$\phi_\bE:\Gamma_{\!\bE}\to\Gamma_{\!\bE'}$,
$\phi_A:\Gamma_{\!A}\to\Gamma_{\!A'}$ fitting into a commutative
diagram :
$$
\xymatrix{
\Gamma_{\!\bE\circ} \ar[d]_{\gamma_{A\circ}}
\ar@/^2.7pc/[rrr]^-{\phi_{\bE\circ}} &
\bE \ar[l]_-{v_\bE} \ar[r]^-{\bE(f)} \ar[d]_{\bar u_A} &
\bE' \ar[r]^{v_{\bE'}} \ar[d]^{\bar u_{A'}} &
\Gamma_{\!\bE'\circ} \ar[d]^{\gamma_{A'\circ}} \\
\Gamma_{\!A\circ} \ar@/_2.7pc/[rrr]_-{\phi_{A\circ}} &
A \ar[l]_-{v_A} \ar[r]^-f & A' \ar[r]^-{v_{A'}} &
\Gamma_{\!A'\circ}
}$$
where $\gamma_A$ and $\gamma_{A'}$ are the group isomorphisms
provided by proposition \ref{prop_A-is-valring-is-E}(ii.b) :
details left to the reader.
\end{remark}

\begin{theorem}\label{th_Scholze-tilt}
Let $A$ be any perfectoid ring, and set $\bE:=\bE(A)$. We have :
\begin{enumerate}
\item
The induced map $\Cont(\bar u_A)$ of proposition
{\em\ref{prop_Cont-u_A}(ii)} is a homeomorphism.
\item
For every rational subset $R$ of\/ $\Cont(\bE)$, the
preimage $\Cont(\bar u_A)^{-1}R$ is a rational subset
of\/ $\Cont(A)$ (see \eqref{subsec_rational-Cont}).
\item
$\Cont(\bar u_A)$ restricts to bijections
$$
\Cont(A)_\mathrm{a}\isom\Cont(\bE)_\mathrm{a}
\qquad\text{and}\qquad
\Cont(A)_\mathrm{na}\isom\Cont(\bE)_\mathrm{na}.
$$
\item
Let $w\in\Cont(A)$ be any valuation, and set
$v:=\Cont(\bar u_A)(w)$. Then $\bar u_A$ induces group
isomorphisms
$$
\gamma_w:\Gamma_{\!v}\isom\Gamma_{\!w}
\qquad\text{and}\qquad
c\Gamma_{\!v}\isom c\Gamma_{\!w}.
$$
\end{enumerate}
\end{theorem}
\begin{proof}(ii): If $A$ is perfectoid, both $\Cont(A)$
and $\Cont(\bE)$ are spectral spaces, and their rational
subsets form a basis of constructible open subsets.
Thus, let $(e_i~|~i=0,\dots,n)$ be any finite system of
elements of $\bE$ that generate an open ideal, and set
$R:=R_\bE(\frac{e_1}{e_0},\dots,\frac{e_n}{e_0})\cap\Cont(\bE)$.
Let also $a_i:=\bar u_A(e_i)$ for every $i=0,\dots,n$. Then
$$
\Cont(\bar u_A)^{-1}R=
R_A\Bigl(\frac{a_1}{a_0},\dots,\frac{a_n}{a_0}\Bigr)\cap\Cont(A).
$$
However, the system $(a_i~|~i=0,\dots,n)$ generates an open ideal
of $A$ (corollary \ref{cor_taut-two}(ii)), whence the assertion.

(iii): Let $\underline\alpha:=(\alpha_n~|~n\in\N)$ be any
distinguished element of $\Ker\,u_A$; since $\alpha_0$ is
topologically nilpotent, it lies in every open prime ideal
of $\bE(A/pA)$. Likewise, $pA$ is contained in every open prime
ideal of $A$, and due to lemmata \ref{lem_was-third-cond}(i,iii)
and \ref{lem_drop-conditions}(iv), the map $\bar u_A$ induces an
isomorphism
$$
\omega:\bE/\alpha_0\bE\isom A/pA
$$
of topological rings, so that $\Cont(\bar u_A)$ agrees with
$\Cont(\omega)$ on the closed subset $\Cont(\bE/\alpha_0\bE)$
of $\Cont(\bE)$, and induces a homeomorphism
$$
\Cont(\bE/\alpha_0\bE)\isom\Cont(A/pA).
$$
It is also clear that $\Cont(\omega)$ maps
$\Cont(\bE/\alpha_0\bE)_\mathrm{na}=\Cont(\bE)_\mathrm{na}$
homeomorphically onto
$\Cont(A/pA)_\mathrm{na}=\Cont(A)_\mathrm{na}$, so we
need only to check that the restriction
$$
C(A):=\Cont(A)\cap R_A\Bigl(\frac{p}{p}\Bigr)\to C(\bE):=
\Cont(\bE)\cap R_\bE\Bigl(\frac{\alpha_0}{\alpha_0}\Bigr)
$$
of $\Cont(\bar u_A)$ is a bijection. Set $C^+(A):=C(A)\cap\Spv^+A$
and $C^+(\bE):=C(\bE)\cap\Spv^+\bE$. Clearly $\Cont(\bar u_A)$
restricts to a mapping $C^+(A)\to C^+(\bE)$, and we shall
first exhibit an inverse mapping
$$
\sigma:C^+(\bE)\to C^+(A)
$$
for $\Cont(\bar u_A)_{|C^+(A)}$. To this aim, let $v\in C^+(\bE)$
be any valuation; then $v$ factors through a ring homomorphism
$\pi:\bE\to V:=\kappa(v)^+$ and the residual valuation
$\bar v:V\to\Gamma_{\!v\circ}$ of $v$.
Since $\bE$ is a perfect $\F_p$-algebra, the field
$\kappa(v)$ is perfect, hence $V$ is a perfect
$\F_p$-algebra as well : indeed, it is clear that
$x^p\in V$ if and only if $x\in V$, for every
$x\in\kappa(v)$. Let us endow $V$ with its valuation
topology $\cT_V$, so that $\bar v$ is continuous for
the topology $\cT_V$ (remark \ref{rem_Cont-A}(i)).

\begin{claim}\label{cl_pi-is-continuous}
The completion $(V^\wedge,\cT^\wedge_V)$ of
$(V,\cT_V)$ is a perfectoid valuation ring.
\end{claim}
\begin{pfclaim} We know that $V^\wedge$ is a valuation ring,
by proposition \ref{prop_stays-valuation}(iii), and its
topology $\cT^\wedge_V$ is adic with a principal ideal of
adic definition, since the same holds for $\cT_V$
(lemma \ref{lem_Cont-A}(ii)).
Moreover, it is easily seen that the Frobenius endomorphism
$\Phi_V$ of $V$ is a homeomorphism for the topology $\cT_V$,
hence it induces a homeomorphism $\Phi^\wedge_V$ on $V^\wedge$;
but clearly $\Phi^\wedge_V$ is just the Frobenius endomorphism
of $V^\wedge$, so the latter is a perfect topological
$\F_p$-algebra, and the assertion follows (see example
\ref{ex_perfectoid}(i)).
\end{pfclaim}

Let $\pi^\wedge:\bE\to V^\wedge$ be the composition of $\pi$
with the completion map $V\to V^\wedge$; since $\pi^\wedge$
is a continuous map of perfectoid $\F_p$-algebras (remark
\ref{rem_Cont-A}(i)), the image of $\underline\alpha$
under the map $W(\pi^\wedge)$ is still distinguished
in $W(V^\wedge)$, so claim \ref{cl_pi-is-continuous} and
proposition \ref{prop_A-is-valring-is-E}(i,ii) say that
$$
A_V:=W(V^\wedge)\otimes_{W(\bE)}A
$$
is a perfectoid valuation ring whose valuation $v_{A_V}$,
is continuous, and the map $\bar u_{A_V}$ induces
an isomorphism of ordered abelian groups from the value
group of $V^\wedge$ onto that of $A_V$
$$
\gamma_{V^\wedge}:\Gamma_{\!V^\wedge}\isom\Gamma_{\!A_V}
\qquad\text{such that}\qquad
v_{A_V}\circ\bar u_{A_V}=\gamma_{V^\wedge\circ}\circ\bar v{}^\wedge
$$
where $\bar v{}^\wedge$ is the valuation of $V^\wedge$.
Moreover, we get a continuous map of perfectoid rings
$$
\phi:=W(\pi^\wedge)\otimes_{W(\bE)}A:A\to A_V
$$
and $w:=v_{A_V}\circ\phi$ is a continuous valuation
of $A$ such that
$$
\Cont(\bar u_A)(w)=w\circ\bar u_A=
v_{A_V}\circ\bar u_{V^\wedge}\circ\pi^\wedge=
\gamma_{V^\wedge\circ}\circ\bar v{}^{\wedge}\circ\pi^\wedge=
\gamma_{V^\wedge\circ}\circ v.
$$
Since $v(\alpha_0)\neq 0$, we have $w(p)\neq 0$ as well,
so we obtain the sought map $\sigma$, by setting
$\sigma(v):=w$, and we already see that $\sigma$ is a
right inverse for $\Cont(\bar u_A)_{|C^+(A)}$. We notice
that since $v(\bE)\cap\Gamma_{\!v}$ generates the group
$\Gamma_{\!v}$, the subset $w\circ\bar u_A(\bE)\cap\Gamma_{\!A_V}$
generates the group $\Gamma_{\!A_V}$, so that $\bar u_A$
induces an isomorphism
\set\begin{equation}\label{eq_almost-iv}
\Gamma_{\!v}\isom\Gamma_{\!\sigma(v)}
\qquad
\text{for every $v\in C^+(\bE)$}.
\end{equation}
Next, we check that $\sigma$ is a left inverse as well.
To this aim, let $v'\in C^+(A)$ be any valuation; then
$v'$ factors through a ring homomorphism
$\pi':A\to V':=\kappa(v')^+$ and the residual valuation
$\bar v{}':V'\to\Gamma_{\!v'\circ}$. By lemma
\ref{lem_Cont-A}(ii), the map $\pi'$ is continuous for
the valuation topology $\cT_{V'}$ of $V'$, and $\cT_{V'}$
is adic with a principal ideal of adic definition.
Since $pV'\neq 0$ and $p$ is topologically nilpotent in
$A$, the ideal $pV'$ is open and topologically nilpotent
in $V'$, so $\cT_{V'}$ agrees with the $p$-adic topology on
$V'$. Let $(V'^\wedge,\cT^\wedge_{V'})$ be the completion of
$(V',\cT_{V'})$; then $V'^\wedge$ is a valuation ring whose
valuation $\bar v{}'^\wedge:V'^\wedge\to\Gamma_{\!v'\circ}$ extends
$\bar v{}'$ (proposition \ref{prop_stays-valuation}(iii,v)),
and we let $\pi'^\wedge:A\to V'^\wedge$ be the composition
of $\pi'$ with the completion map $V'\to V'^\wedge$.
Then $V'_\bE:=\bE(V'^\wedge)$ is a valuation ring whose
valuation is
$$
\bar v{}'_\bE:=\bar v{}'^\wedge\circ\bar u_{V'^\wedge}:
V'_\bE\to\Gamma_{\!v'\circ}
$$
(lemma \ref{lem_E-is-val.ring}) and $\bE(\pi'^\wedge):\bE\to V'_\bE$
is a continuous ring homomorphism. Since
$\pi'^\wedge\circ\bar u_A=\bar u_{V'^\wedge}\circ\bE(\pi'^\wedge)$,
we deduce that
\set\begin{equation}\label{eq_equality-of-classes}
v'_\bE:=\bar v{}'_\bE\circ\bE(\pi'^\wedge)=\Cont(\bar u_A)(v')
\qquad
\text{in $\Cont(\bE)$}.
\end{equation}

\begin{claim}\label{cl_V-prime-is-perfectoid}
(i)\ \ 
The topology of $V'_\bE$ agrees with the valuation topology
$\cT_{V'_\bE}$.
\begin{itemize}
\item[(ii)]
$(V'_\bE,\cT_{V'_\bE})$ is a perfectoid ring.
\end{itemize}
\end{claim}
\begin{pfclaim}(i): Since the quotient topology induced by
$\cT_{V'}$ on $V'/pV'$ is discrete, the topology of $V'_\bE$
is linear, complete and separated (remark
\ref{rem_topology-of-E}(i)),
so it suffices to check that this topology is not discrete
(proposition \ref{prop_stays-valuation}(i)). To this aim,
since $\alpha_0$ is topologically nilpotent in $\bE$, it
suffices to check that $\nu:=\bE(\pi'^\wedge)(\alpha_0)\neq 0$.
However :
$$
\bar u_{V'^\wedge}(\nu)=
\pi'^\wedge\circ\bar u_A(\alpha_0)=\pi'^\wedge(pu)
\qquad
\text{for some $u\in A^\times$}
$$
(lemma \ref{lem_was-third-cond}(iii)), and since $v'(p)\neq 0$,
we must have $\pi'^\wedge(p)\neq 0$, whence the claim.

(ii): We have just seen that $\nu$ is a non-zero topologically
nilpotent element of $V'_\bE$, so $\cT_{V'_\bE}$ agrees with
the $\nu$-adic topology on $V'_\bE$; then the assertion
follows from example \ref{ex_perfectoid}(i).
\end{pfclaim}

Arguing as in the foregoing we see that the image of
$\underline\alpha$ in $W(V'_\bE)$ is still distinguished, so
$$
A_{V'}:=W(V'_\bE)\otimes_{W(\bE)}A
$$
is a perfectoid valuation ring whose value group is naturally
identified with the value group $\Gamma_{\!v'_\bE}$ of $v'_\bE$
(proposition \ref{prop_A-is-valring-is-E}(i,ii.b) and claim
\ref{cl_V-prime-is-perfectoid}(ii)); by the same token,
we get a natural isomorphism
$$
\beta:\bE(A_{V'})\isom V'_\bE
\qquad\text{such that}\qquad
\pi_{V'_\bE}\circ W(\beta)=u_{A_{V'}}
$$
where $\pi_{V'_\bE}:W(V'_\bE)\to A_{V'}$ is the projection
(proposition \ref{prop_special-case}(iv)).
Moreover, in light of the commutative diagram
$$
\xymatrix{ W(\bE) \ar[r]^-{u_A} \ar[d]_{W(\bE(\pi'^\wedge))} &
A \ar[d]^{\pi'^\wedge} \\
W(V'_\bE) \ar[r]^-{u_{V'^\wedge}} & V'^\wedge
}$$
we get a unique homomorphism of $A$-algebras
$$
\psi:A_{V'}\to V'^\wedge
\qquad\text{such that}\qquad
u_{V'^\wedge}=\psi\circ\pi_{V'_\bE}.
$$
There follows a diagram :
$$
\xymatrix{ A_{V'} \ar[d]_{v'_A} \ar@/^2.7pc/[rrrr]^-\psi & &
V'_\bE \ar[ll]_-{\bar u_{A_{V'}}\circ\beta^{-1}}
\ar[rr]^-{\bar u_{V'^\wedge}} \ar[d]^{\bar v{}'_\bE} & &
V'^\wedge \ar[d]^{\bar v{}'^\wedge} \\
\Gamma_{\!v'_\bE\circ} \rrdouble & &
\Gamma_{\!v'_\bE\circ} \ar[rr]^-j & & \Gamma_{\!v'\circ}
}$$
whose two square subdiagrams commute, where $v'_A$ is the
valuation of $A_{V'}$ and $j$ the inclusion map of
$\Gamma_{\!v'_\bE\circ}$ into $\Gamma_{\!v'\circ}$. Let also
$\phi':A\to A_{V'}$ be the structure map of the $A$-algebra
$A_{V'}$.

\begin{claim}\label{cl_twice-che-fatica}
(i)\ \ The top subdiagram commutes as well.
\begin{enumerate}
\addenu
\item
The external subdiagram also commutes, {\em i.e.}
$\bar v{}'^\wedge\circ\psi=j\circ v'_A$.
\item
The topology of $A_{V'}$ agrees with its valuation topology.
\item
$\psi$ is an isomorphism of topological rings and
$\Gamma_{\!v'_\bE}=\Gamma_{\!v'}$.
\end{enumerate}
\end{claim}
\begin{pfclaim}(i): Indeed, taking into account lemma
\ref{lem_drop-conditions}(iv) we may compute :
$$
\begin{aligned}
\psi\circ\bar u_{A_{V'}}=&\,
\psi\circ u_{A_{V'}}\circ\tau_{\bE(A_{V'})} \\
=&\,
\psi\circ\pi_{V'_\bE}\circ W(\beta)\circ\tau_{\bE(A_{V'})} \\
=&\,u_{V'}\circ W(\beta)\circ\tau_{\bE(A_{V'})} \\
=&\,u_{V'}\circ\tau_{V'_\bE}\circ\beta \\
=&\,\bar u_{V'}\circ\beta.
\end{aligned}
$$

(ii): Let $x\in A_{V'}$ be any element; since the left
square subdiagram commutes, we may find $y\in V'_\bE$
such that $\bar u_{A_{V'}}\circ\beta^{-1}(y)\cdot u=x$ for
some $u\in A_{V'}^\times$, and therefore
$\psi(x)=\bar u_{V'^\wedge}(y)\cdot\psi(u)$. Obviously,
$\psi(u)$ is invertible in $V'^\wedge$, so that
$\bar v{}'^\wedge\circ\psi(x)=
\bar v{}'^\wedge\circ\bar u_{V'^\wedge}(y)=
j\circ\bar v{}'_\bE(y)=j\circ v'_A(x)$, as stated.

(iii): In light of proposition \ref{prop_stays-valuation}(i),
it suffices to check that the topology of $A_{V'}$ is not
discrete. However, notice that by construction $\phi'$ is
a continuous map; thus, if the topology of $A_{V'}$ were
discrete, $\Ker\,\phi'$ would be an open ideal, and then
the same would hold for $\Ker\,\pi'^\wedge$. But this is
absurd, since $v'$ is analytic.

(iv): As an immediate consequence of (ii), we already see that
$\psi$ is injective. Then it is clear that $\phi'$ factors
uniquely through $\pi'$ and an injective map $V'\to A_{V'}$
whose composition with $\psi$ equals the completion map.
Thus, the image of $\psi$ is a dense subring of $V'^\wedge$;
combining with (ii), we deduce already the stated equality
of value groups. Then by (iii) the topology of $A_{V'}$ agrees
with the topology induced by $V'^\wedge$ via $\psi$. Since
$A_{V'}$ is complete and separated, theorem
\ref{th_complete-top-grps}(iii) implies that $\psi$ is an
isomorphism.
\end{pfclaim}

From claim \ref{cl_twice-che-fatica}(ii) we get the following
equalities in $\Cont(A)$ :
$$
\begin{aligned}
\sigma\circ\Cont(\bar u_A)(v')=\,&
\sigma(\bar v{}'^\wedge\circ\bar u_{V'^\wedge}\circ\bE(\pi'^\wedge_\bE)) \\
=\,& \sigma(\bar v{}'_\bE\circ\bE(\pi'^\wedge)) \\
=\,& v'_A\circ\phi' \\
=\,& \bar v{}'^\wedge\circ\psi\circ\phi' \\
=\,& \bar v{}'^\wedge\circ\pi'^\wedge \\
=\,& v'
\end{aligned}
$$
as required. Let $\kappa(v'_\bE)^{\wedge+}$ be the
completion of $\kappa(v'_\bE)^+$ for its valuation
topology, $\pi'_\bE:\bE\to\kappa(v'_\bE)^+$ the natural
ring homomorphism (whose composition with the valuation
$\kappa(v'_\bE)^+\to\Gamma_{\!v'\circ}$ yields $v'_\bE$),
and $\pi'^\wedge_\bE:\bE\to\kappa(v'_\bE)^{\wedge+}$ the
completion of $\pi'_\bE$; we notice as well :

\begin{claim}\label{cl_compaq-sold}
There exists a ring isomorphism
$$
\omega^\wedge_\bE:\kappa(v'_\bE)^{\wedge+}\isom V'_\bE
\qquad\text{such that}\qquad
\omega^\wedge_\bE\circ\pi'^\wedge_\bE=\bE(\pi'^\wedge).
$$
\end{claim}
\begin{pfclaim} The map $\bE(\pi'^\wedge)$ factors through
$\pi'_\bE$ and a unique injective ring homomorphism
$\omega_\bE:\kappa(v'_\bE)^+\to V'_\bE$ such that
$\bar v{}'_\bE\circ\omega_\bE:\kappa(v'_\bE)^+\to\Gamma_{\!v'\circ}$
is the residual valuation. Thus, the valuation topology
of $\kappa(v'_\bE)^+$ agrees with the topology induced from
$V'_\bE$ via $\omega_\bE$, so the latter extends to a continuous
ring homomorphism
$\omega^\wedge_\bE:\kappa(v'_\bE)^{\wedge+}\to V'_\bE$ such that
$\bar v{}'_\bE\circ\omega^\wedge_\bE$ is the valuation of
$\kappa(v'_\bE)^{\wedge+}$. The map $\omega^\wedge_\bE$ is
injective (proposition \ref{prop_replaces-Mat-Th.8.1}(i)),
and clearly $\omega^\wedge_\bE\circ\pi'^\wedge_\bE=\bE(\pi'^\wedge)$,
so it remains only to check that $\omega^\wedge_\bE$ is bijective.
To this aim, set $V'':=W(\kappa(v'_\bE)^{\wedge+})\otimes_{W(\bE)}A$;
by the foregoing we know already that $V''$ is a perfectoid
valuation ring with value group $\Gamma_{\!v'}$, and it suffices
to show that
$$
\omega_A^\wedge:=W(\omega^\wedge_\bE)\otimes_{W(\bE)}A:V''\to A_{V'}
$$
is an isomorphism. However, set
$\pi'^\wedge_A:=W(\pi'^\wedge_\bE)\otimes_{W(\bE)}A$; by remark
\ref{rem_A-is-valring-is-E}, we know already that
$\omega^\wedge_A$ is injective, and moreover we have a
commutative diagram
$$
\xymatrix{ A \ar[r]^-{\pi'^\wedge_A} \ar[rd]_{\phi'} &
V'' \ar[r]^-{v''} \ar[d]^{\omega^\wedge_A} &
\Gamma_{\!v'\circ} \ddouble \\
& A_{V'} \ar[r]^-{v'_A} & \Gamma_{\!v'\circ}
}$$
where $v''$ is the valuation of $V''$. It follows that the
valuation topology of $V''$ is induced by the valuation
topology of $A_{V'}$ via $\omega^\wedge_A$. Furthermore, we
deduce that $v''\circ\pi'^\wedge_A=v'_A\circ\phi'=v'$, so
$\pi'^\wedge_A$ factors through $\pi':A\to V'$ and an
injective map $V'\to V''$ whose composition with
$\omega^\wedge_A$ has dense image, by virtue of claim
\ref{cl_twice-che-fatica}(iv). Since $V''$ is complete, it
follows that $\omega^\wedge_A$ is surjective, whence the claim.
\end{pfclaim}

Next, let $A^+$ be the smallest subring of integral
elements of $A$ (remark
\ref{rem_shouldbe-quasi-affinoid}(iv)); by theorem
\ref{th_int-subrings-perfectoid}(ii,iii), the ring
$A^+$ is perfectoid, and $\bE(A^+)$ is naturally
identified with the smallest ring $\bE^+$ of integral
elements of $\bE$. Moreover, we get a commutative
diagram
$$
\xymatrix{ \Cont(A) \ar[rr]^-{\Cont(\bar u_A)}
\ar[d]_{\Cont(j_A)} & & \Cont(\bE) \ar[d]^-{\Cont(j_\bE)} \\
\Cont(A^+) \ar[rr]^-{\Cont(\bar u_{A^+})} & & \Cont(\bE^+)
}$$
where $j_A:A^+\to A$ and $j_\bE:\bE^+\to\bE$ are the
inclusion maps. Since $A^+$ is the integral closure
in $A$ of the $\Z$-subalgebra generated by $A^{\circ\circ}$,
it is easily seen that $\Cont(A^+)=\Cont^+(A^+)$, and
likewise for $\Cont(\bE^+)$. Therefore, $C(A^+)=C^+(A^+)$
and $C(\bE^+)=C^+(\bE^+)$, and the foregoing shows that
$\Cont(\bar u_{A^+})$ restricts to a bijection
$$
C(\bE^+)\isom C(A^+).
$$
Lastly, proposition \ref{prop_Cont-open-subring}(ii) says
that $\Cont(j_A)$ restricts to a bijection from the analytic
valuations of $A$ to those of $A^+$, and likewise for
$\Cont(j_\bE)$, whence (iii).

(iv): Let $\pi_A:A\to A/pA$ and $\pi_\bE:\bE\to\bE/\alpha_0\bE$
be the projections; if $w(p)=0$, then $w=\Cont(\pi_A)(w')$
for a unique $w'\in\Cont(A/pA)$, and $v=\Cont(\pi_\bE)(v')$,
where $v':=\Cont(\omega)(w')$. Since $\pi_A$ and $\pi_\bE$ are
surjective, obviously $\Gamma_{\!w}=\Gamma_{\!w'}$ and
$\Gamma_{\!v}=\Gamma_{\!v'}$, and likewise for the characteristic
subgroups; since $\omega$ is a ring isomorphism, it induces
a group isomorphism $\Gamma_{\!v'}\isom\Gamma_{\!w'}$ that
identifies the respective characteristic subgroups, whence
the assertion, in this case. Next, suppose that $w\in C(A)$,
and set $w':=\Cont(j_A)(w)$, $v':=\Cont(\bar u_{A^+})(w')$;
then $w'\in C^+(A^+)$, and from (iii) and \eqref{eq_almost-iv}
we see that $\bar u_{A^+}$ induces an isomorphism
$\gamma_{w'}:\Gamma_{\!w'}\isom\Gamma_{\!v'}$. Moreover, we
have a commutative diagram of groups
$$
\xymatrix{ \Gamma_{\!w'} \ar[r]^-{\gamma_{w'}} \ar[d]_{i_w} &
\Gamma_{\!v'} \ar[d]^{i_v} \\
\Gamma_{\!w} \ar[r]^-{\gamma_w} & \Gamma_{\!v}
}$$
where $i_w$ and $i_v$ are the injective maps induced
by $j_A$ and $j_\bE$. In order to prove that $\gamma_{\!w}$
is an isomorphism, it suffices therefore to check that
$i_w$ and $i_v$ are surjective. To this aim, it suffices
to show that $S_{A^+}:=w(A^+)\cap\Gamma_{\!w}$ generates
$\Gamma_{\!w}$ and $S_{\bE^+}:=v(A^+)\cap\Gamma_{\!v}$
generates $\Gamma_{\!v}$. However, notice that $pA\subset A^+$
and $\alpha_0\bE\subset\bE^+$; since $w(p)\neq 0$ and
$v(\alpha_0)\neq 0$, it follows easily that $S_{A^+}$ and
$S_A:=w(A)\cap\Gamma_{\!w}$ generate the same subgroup of
$\Gamma_{\!w}$, and likewise for $S_{\bE^+}$. To conclude,
it suffices to remark that $S_A$ generates $\Gamma_{\!w}$,
and likewise for $\Gamma_{\!v}$.

Next, obviously $\gamma_w$ restricts to an injective
group homomorphism $c\Gamma_{\!v}\to c\Gamma_{\!w}$.
To show that this map is surjective, consider any
$a\in A$ such that $w(a)>1$; since $\bar u_{A/pA}$ is
surjective (lemma \ref{lem_was-third-cond}(i)), we may
write $a=\bar u_A(\beta)+pb$ for some $\beta\in\bE$
and $b\in A$, and since $pb\in A^{\circ\circ}$, we have
$w(pb)<1$, so that $w(a)=v(\beta)$, whence the
contention.

(i): We argue as in the proof of proposition
\ref{prop_Cont-of-complete} : by (ii) and (iii) we
know already that $\Cont(\bar u_A)$ is a spectral and
bijective map of spectral spaces, so it suffices to
check that $\Cont(\bar u_A)$ is a closed map, and by
proposition \ref{prop_closed-under-spec}(i) and corollary
\ref{cor-pro-constr}(i) we are further reduced to showing
that $\Cont(\bar u_A)$ is specializing. Hence, let
$w\in\Cont(A)$ be any valuation, and $v'\in\Cont(A)$
any specialization of $v:=\Cont(\bar u_A)(w)$; we may
find a valuation $v''$ of $\bE$ that is both a secondary
specialization of $v$ and a primary generization of $v'$,
and $v''\in\Cont(\bE)$, by remark \ref{rem_Cont-A}(iv).
By proposition \ref{prop_Cont-u_A}(v), there exists a
secondary specialization $w''$ of $w$ such that
$\Cont(\bar u_A)(w'')=v''$. By (iv) we know that $w''$
and $v''$ have the same value groups and characteristic
subgroups, so we may find a primary specialization $w'$
of $w$ such that $\Cont(\bar u_A)(w')=v'$, whence the
assertion.
\end{proof}

\sset\subsubsection{}\label{subsec_what-name}
Let $\underline X$ be any perfectoid quasi-affinoid
scheme, and with the notation of \eqref{subsec_Gamma-circ}
and \eqref{subsec_upgrade-E}, let us set
$$
\underline A:=(A^\circ_X,A^+_X,X):=\sGamma^\circ(\underline X)
\qquad\text{and}\qquad
\underline\bE:=(\bE_X^\circ,\bE_X^+,X_\bE):=\bE(\underline A).
$$
Especially, $\bE_X^\circ=\bE(A^\circ_X)$ and $\bE_X^+=\bE(A^+_X)$;
it is then clear that $\Cont(\bar u_{A^\circ_X})$ maps the subset
$\Spa\,(A^\circ_X,A^+_X)$ into $\Spa\,(\bE_X^\circ,\bE_X^+)$.
Conversely, let $v:A^\circ_X\to\Gamma_v$ be any continuous
valuation such that
$v\circ\bar u_{A^\circ_X}\in\Spa\,(\bE_X^\circ,\bE_X^+)$, and
$\underline\alpha$ any distinguished element of
$\Ker\,u_{A^\circ_X}$; for every $a\in A^+_X$ there exists
$b\in A^+_X$ and $e\in\bE^+_X$ such that
$a=\bar u_{A^\circ_X}(e)+\alpha_0b$, and since
$\alpha_0b\in A^{\circ\circ}_X$, we deduce that $v(a)\leq 1$,
so $v\in\Spa\,(A^\circ_X,A^+_X)$. Combining with proposition
\ref{prop_Cont-u_A}(iii), we deduce that
$\Cont(\bar u_{A^\circ_X})$ restricts to a homeomorphism
$$
\Spa\,(\bar u_{\underline A}):
\Spa\,\underline A\isom\Spa\,\underline\bE
$$
and taking into account the natural identifications
$\Spa\,\underline A\isom\Spa\,\underline X$ and
$\Spa\,\underline\bE\isom\Spa\,\bE(\underline X)$, we get
as well an induced homeomorphism 
$$
\Spa(\bar u_{\underline X}):
\Spa\,\underline X\isom\Spa\,\bE(\underline X).
$$
Moreover, it is easily seen that every morphism
$f:\underline Y\to\underline X$ of perfectoid quasi-affinoid
schemes yields a commutative diagram
\set\begin{equation}\label{eq_naturality-bE}
{\diagram
\Spa\,\underline Y \ar[rr]^-{\Spa(\bar u_{\underline Y})}
\ar[d]_{\Spa\,f}
& & \Spa\,\bE(\underline Y) \ar[d]^{\Spa\,\bE(f)} \\
\Spa\,\underline X \ar[rr]^-{\Spa(\bar u_{\underline X})} & &
\Spa\,\bE(\underline X).
\enddiagram}
\end{equation}

\begin{corollary}\label{cor_vals-are-perfectoid}
In the situation of \eqref{subsec_what-name}, let
$x\in\Spa\,\underline X$ be any element and set
$y:=\Spa\,(\bar u_{\underline X})(x)$. Denote also by
$\pi^{\wedge+}_x:A^+_X\to\kappa(x)^{\wedge+}$ and
$\pi^{\wedge+}_y:\bE^+_X\to\kappa(y)^{\wedge+}$ the natural
maps (notation of \eqref{subsec_valuation-on-stalks}),
and endow $\kappa(x)^{\wedge+}$ and $\kappa(y)^{\wedge+}$
with the unique ring topologies $\cT^\wedge_x$ and
$\cT^\wedge_y$ such that $\pi^{\wedge+}_x$ and $\pi^{\wedge+}_y$
are adic ring homomorphisms. We have :
\begin{enumerate}
\item
The topological rings $(\kappa(x)^{\wedge+},\cT^\wedge_x)$ and
$(\kappa(y)^{\wedge+},\cT^\wedge_y)$ are perfectoid valuation
rings.
\item
There exists a unique isomorphism of topological rings :
$$
\omega^{\wedge+}:\bE(\kappa(x)^{\wedge+})\isom\kappa(y)^{\wedge+}
\qquad\text{such that}\qquad
\omega^{\wedge+}\circ\bE(\pi^{\wedge+}_x)=\pi^{\wedge+}_y.
$$
\end{enumerate}
\end{corollary}
\begin{proof} Suppose first that $p\in\Ker\,\pi^+_x$; in this
case $x$ corresponds to a continuous ring homomorphism
$\bar\pi_x:A^\circ_X/pA^\circ_X\to\kappa(x)$ and $y$ corresponds
to the continuous ring homomorphism
$\bar\pi_y:=\bar\pi_x\circ\phi:
\bE^\circ_X/\alpha_0\bE^\circ_X\to\kappa(x)$, where $\alpha_0$ is
as in \eqref{subsec_what-name} and
$\phi:\bE^\circ_X/\alpha_0\bE^\circ_X\isom A^\circ_X/pA^\circ_X$ is
the isomorphism of topological rings induced by $\bar u_{A^\circ_X}$.
There follows a unique isomorphism of valued fields
$\omega:\kappa(y)\isom\kappa(x)$ fitting into the commutative
diagram
$$
\xymatrix{
\bE^\circ_X/\alpha_0\bE^\circ_X \ar[r]^-\phi \ar[d]_{\bar\pi_y} &
A^\circ_X/pA^\circ_X \ar[d]^{\bar\pi_x} \\
\kappa(y) \ar[r]^-\omega & \kappa(x).
}$$
Let us now endow $\kappa(x)^+$ and $\kappa(y)^+$ with the unique
ring topologies such that the natural maps
$\pi^+_x:A^+_X\to\kappa(x)^+$ and $\pi^+_y:\bE^+_X\to\kappa(y)^+$
are adic; then clearly $\omega$ restricts to an isomorphism
$\omega^+:\kappa(y)^+\isom\kappa(x)^+$ of topological rings,
and assertion (ii) follows straightforwardly in this case.
Moreover, since $\bE$ is a perfect $\F_p$-algebra, the same
holds for its quotient $\bE/\Ker\,\pi_y$, and then also for
$\kappa(y)=\Frac\,(\bE/\Ker\,\pi_y)$, and for $\kappa(y)^\wedge$
as well (example \ref{ex_discrete-Witt}(ii)). By remark
\ref{rem_topology-of-E}(v), it follows easily that
$\kappa(y)^{\wedge+}$ is a perfect adic and f-adic topological
ring (details left to the reader), and hence it is perfectoid
(example \ref{ex_perfectoid}(i)); the same then holds for
$\kappa(x)^{\wedge+}$.

Next, if $\pi^+_x(p)\neq 0$, the point $x$ is analytic,
and the same then holds for $y$. In this case, it is easily
seen that the topologies $\cT^\wedge_x$ and $\cT^\wedge_y$ are
separated and not discrete, so they coincide with the respective
valuation topologies (proposition \ref{prop_stays-valuation}(i)).
Recall that the inclusion map $i:A^+_X\to A^\circ_X$ induces a
homeomorphism $\Cont(i)_\mathrm{a}:
\Cont(A^\circ_X)_\mathrm{a}\isom\Cont(A^+_X)_\mathrm{a}$ and
likewise for $\Cont(\bE^\circ_X)$ (proposition
\ref{prop_Cont-open-subring}(ii)); moreover, if
$x':=\Cont(i)(x)$, we have a unique isomorphism of valued
fields $\kappa(x')\isom\kappa(x)$ fitting into the commutative
diagram
$$
\xymatrix{ A^+_X \ar[rr]^-i \ar[d]_{\pi^+_{x'}} & &
A^\circ_X \ar[d]^{\pi^\circ_x} \\
\kappa(x')^+ \ar@{^{(}->}[r] & \kappa(x') \ar[r]^-\sim & \kappa(x)
}$$
where $\pi^\circ_x$ and $\pi^+_{x'}$ are the natural maps,
and likewise for $y$ and its image $y'\in\Cont(\bE^+_X)$.
Hence we may regard $x$ and $y$ as elements of
$\Cont(A^+_X)$ and respectively $\Cont(\bE^+_X)$, and
in this case assertions (i) and (ii) follow from claims
\ref{cl_twice-che-fatica}(iv) and \ref{cl_compaq-sold}. 
\end{proof}

\subsection{Graded perfectoid rings}
This section is mainly dedicated to the construction and
investigation of certain {\em angular Rees algebras} that
shall be useful in the study of blowing up morphisms of
perfectoid formal schemes.

\begin{proposition}\label{prop_graded-perfectoid}
Let $(\Gamma,+,0)$ be a monoid, $(A,\underline B)$ a
$\Gamma$-graded structure on the topological ring
$(A,\cT)$, and $i_0:\gr_0B\to A$ the inclusion map.
Suppose that $A$ is a P-ring. Then we have :
\begin{enumerate}
\item
If the $p$-Frobenius endomorphism $\bp_\Gamma$ of\/
$\Gamma$ is injective, there exists a finitely generated
graded ideal $J$ of $B$ such that $JA$ is an ideal of
definition of $A$.
\item
Suppose that $i_0$ is c-adic and $\bp_\Gamma$ is injective.
Then $\gr_0B$ is a P-ring, and if $A$ is perfectoid, the same
holds for $\gr_0B$.
\item
If $A$ is a perfectoid ring, $\Ker\,u_A$ is generated
by a distinguished element in $\bA(\gr_0B)$.
\end{enumerate}
\end{proposition}
\begin{proof}(i): By lemma \ref{lem_perfectoid}(iv), we may
write $p=b^pu$ for some $u\in A^\times$ and $b\in A$, and
by the same token, $u=x^p+py$ for some $x,y\in A$,
so that $p=a^p+p^2w$, with $a:=bx$ and $w:=u^{-1}y$.
Let now $(a_\gamma~|~\gamma\in\Gamma)$ be the sequence
attached to $a$, as in remark \ref{rem_graded-top-algs}(iii),
and $c$ (resp. $d$) the limit of the Cauchy net
$(c_S~|~S\subset\Gamma)$ (resp. $(d_S~|~S\subset\Gamma)$)
provided by lemma \ref{lem__p-power-Cauchy}(i), where $S$
ranges over the finite subsets of $\Gamma$, and 
$a_S^p=c_S+p\cdot d_S$ for every such $S$ (with $a_S$
defined as in remark \ref{rem_graded-top-algs}(iii)).
In light of proposition \ref{prop_Cauchy}(ii), we deduce
that $a^p=c+pd$, whence $p=c+pd+p^2w$. Moreover, since the
$p$-Frobenius map of $\Gamma$ is injective, we have
$\pi'_{p\gamma}(c_S)=a_\gamma^p$ for every $S\subset\Gamma$
and every $\gamma\in S$; therefore $\pi'_{p\gamma}(c)=a_\gamma^p$
for every $\gamma\in\Gamma$, and especially, $\pi'_0(c)=a_0^p$.
Consequently
$$
p=a_0^p+p\cdot d_0+p^2w_0
$$
where $(d_\gamma~|~\gamma\in\Gamma)$ and
$(w_\gamma~|~\gamma\in\Gamma)$ are the sequences attached
to $d$ and respectively $w$. Next, notice that $a$ is
topologically nilpotent in $A$, hence $a_\gamma$ is
topologically nilpotent in $B$, for every $\gamma\in\Gamma$
(lemma \ref{lem__p-power-Cauchy}(ii)), and then the same
holds for $d_0$. Thus $v:=1-d_0-pw_0$ is invertible in $A$,
with inverse given by the convergent series
$\sum_{n\in\N}(d_0+pw_0)^n$, and since $\gr_0B$ is closed
in $A$, we conclude that $v\in(\gr_0B)^\times$.
Summing up, we arrive at the identity
$$
pv=a_0^p
\qquad
\text{with $v\in(\gr_0B)^\times$ and $a_0\in\gr_0B$}.
$$
Lastly, by proposition \ref{prop_general-graded-top}(i.a)
we know already that $A$ admits an ideal of adic definition
of the form $J'A$, for some graded ideal $J'\subset B$, and
we set $J:=J'+a_0B$; then it is easily seen that $JA$ is an
ideal of definition for the P-ring $A$.

(ii): We remark :

\begin{claim}\label{cl_graded-frob}
In the situation of the proposition, suppose that
$\bp_\Gamma$ is injective. Then we have :
\begin{enumerate}
\item
The Frobenius endomorphism of
$B/pB$ restricts to surjective maps
$$
\gr_\gamma B\otimes_\Z\F_p\to\gr_{p\gamma}B\otimes_\Z\F_p
\qquad
\text{for every $\gamma\in\Gamma$}.
$$
\item
$\gr_\gamma B=0$ for every $\gamma\in\Gamma\setminus\bE(\Gamma)$.
\end{enumerate}
\end{claim}
\begin{pfclaim}(i): Say that $b\in\gr_\gamma B$; by lemma
\ref{lem_perfectoid}(iv) there exists $x\in A$ such that
$x^p-b\in pA$. Let $(x_\gamma~|~\gamma\in\Gamma)$ be the
system attached to $x$ as in remark \ref{rem_graded-top-algs}(iii);
it was remarked in the proof of (i) that the system
$(x_\gamma^p~|~\gamma\in\N)$ is a Cauchy net in $B$, and
if $y\in A$ denotes the limit of this net, then $x^p-y\in pA$,
so that $y-b\in pA$ as well. Consequently,
$$
b=\pi'_\gamma(b)\equiv\pi'_\gamma(y)=
x^p_{p^{-1}\gamma}\pmod{p\cdot\gr_\gamma B}
$$
if $\gamma\in p\Gamma$, and otherwise $b\in p\cdot\gr_\gamma B$.
Assertion (i) follows already. Next, using part (i) of
the claim it is easily seen that, if
$\gamma\in\Gamma\setminus\bE(\Gamma)$, we have
$\gr_\gamma B\otimes_\Z\F_p=0$, whence (ii), since
the $p$-adic topology is separated on $B$.
\end{pfclaim}

Claim \ref{cl_graded-frob} says especially that the Frobenius
endomorphism of $\gr_0B\otimes_\Z\F_p$ is surjective. On the
other hand, we already know by proposition
\ref{prop_general-graded-top}(iii) that the topology of
$\gr_0B$ is adic with a finitely generated ideal of adic
definition $I$. Moreover, the proof of (i) shows that there
exists $a_0\in\gr_0B$ such that $p=a_0^pu$ for some
$u\in(\gr_0B)^\times$; then, it is easily seen that $\gr_0B$
is a P-ring with ideal of definition $I+a_0\cdot\gr_0B$.

Lastly, suppose that $A$ is perfectoid, and to ease notation
set $B_0:=\gr_0B$. By the foregoing, we know already that
$B_0$ is a P-ring, and then $i_0$ is adic, by lemma
\ref{lem_f-adics}(i.c). Pick any ideal of definition $I_0$
of $B_0$, and say that $I_0$ is generated by $k$ elements,
for some $k\in\N$; by remark \ref{rem_p-can-lie-deep}, we
can assume that $p\in I_0^{k(p-1)+1}$, and therefore $I:=I_0A$
is an open ideal of $A$ that fulfills the same condition. Set
$J_0:=I^{(p)}_0$ and $J:=I^{(p)}=J_0A$; it follows that the
Frobenius endomorphism of $A/pA$ induces an isomorphism
$\Phi_I:\gr_I^\bullet A\to\gr_JA$ (propositions
\ref{prop_change-topol}(i) and
\ref{prop_begin-criterion-perf}(ii)).

\begin{claim}\label{cl_get-down-to-B}
The natural maps
$$
\gr_I^nB:=I^n_0B/I^{n+1}_0B\to\gr_I^nA
\qquad
\gr_J^nB:=J^n_0B/J^{n+1}_0B\to\gr_J^nA
$$
are bijective for every $n\in\N$.
\end{claim}
\begin{pfclaim} By claim \ref{cl_exchange-pow-clo} we have
$(I^nB)^c=I^nA$, so that $I^n=I^nB+I^{n+1}$, whence the
surjectivity of the first map. For the injectivity, we need
to show that $I^{n+1}_0B=I^n_0B\cap I^{n+1}$. Hence, let
$x\in I^n_0B\cap I^{n+1}$; then $x=\sum_{i=1}^rx_ra_r$ for
some $r\in\N$ and certain $x_i\in I_0^{n+1}$ and $a_i\in A$
($i=1,\dots,r$).
Say that $x\in\gr_SB:=\bigoplus_{\gamma\in S}\gr_\gamma B$ for
some finite subset $S\subset\Gamma$, and let
$\pi'_S:A\to\gr_SB$ be the induced projection (the sum of
the canonical $\gamma$-projections as in remark
\ref{rem_graded-top-algs}(iii), for $\gamma$ ranging
over the subset $S$). Then
$x=\pi'_S(x)=\sum_{i=0}^rx_i\cdot\pi'_S(a)$, whence the
contention. The same argument applies to the second map.
\end{pfclaim}

From claim \ref{cl_get-down-to-B} we conclude that the
map $\Phi_I:\gr^\bullet_IB\to\gr^\bullet_JB$ is also an
isomorphism. But then, since $\bp_\Gamma$ is injective,
clearly the direct summand $\gr_0\Phi_I:=
\Phi_{I_0}:\gr^\bullet_{I_0}B_0\to\gr^\bullet_{J_0}B_0$ is
still an isomorphism, and therefore $B_0$ is perfectoid,
by theorem \ref{th_criterium-perfect}.

(iii): By proposition \ref{prop_change-topol}(i), we may
assume that $\cT$ is the $p$-adic topology of $A$.  Let
$j_A:A\to B'$ be as in remark \ref{rem_graded-top-algs}(iii);
clearly $p^nB'\cap B_0=p^nB_0$ and
$p^nA\subset j_A^{-1}(p^nB')$, so $B_0\cap p^nA=p^nB_0$
for every $n\in\N$, and therefore the topology of $B_0$
agrees with its $p$-adic topology. It follows that $i_0$
is an adic map, and therefore $B_0$ is even perfectoid,
by (ii); then any distinguished element of $\Ker\,u_{B_0}$
generates $\Ker\,u_A$, by proposition
\ref{prop_cond-for-perfectoid}(ii).
\end{proof}

\sset\subsubsection{}\label{subsec_graded-a_un-cond}
Let $(\Gamma,0,+)$ be a $p$-perfect monoid, $(A,\underline B)$
a perfectoid ring with $\Gamma$-graded structure, and set
$(\bE,\underline B{}_\bE):=\bE(A,\underline B)$. Let also
$\beta_1,\dots,\beta_r$ (resp. $\gamma_1,\dots,\gamma_r$)
be a finite system of elements of $\bE$ (resp. of $\Gamma$),
such that $\beta_i\in\gr_{\gamma_i}\underline B{}_\bE$ for every
$i=1,\dots,r$. Set
$\bff:=(\bar u_A(\beta_1),\dots,\bar u_A(\beta_r))$ and define
the ring $R_{r,0}$ and its ideals $I^{\lceil s\rceil}_{r,0}$ for
every $s\in\R_+$, as in \eqref{subsec_link-with-a-un-etc};
we endow $\underline B$ with the $R_{r,0}$-module structure
induced by the unique ring homomorphism  $u:R_{r,0}\to\underline B$
such that $u(T_i^{1/p^n}):=\bar u_A(\beta^{1/p^n}_i)$ for every
$i=1,\dots,r$ and every $n\in\N$.

\begin{proposition}\label{prop_a-unif-for-perf}
In the situation of \eqref{subsec_graded-a_un-cond},
suppose moreover that the translation map
$\Gamma\to\Gamma\ :\ \gamma\mapsto\gamma+\gamma_i$
is injective for every $i=1,\dots,r$. Then we have :
\begin{enumerate}
\item
$\Tor_i^{R_{r,0}}(R_{r,0}/I^{\lceil s\rceil}_{r,0},\underline B)=0$
for every $i>0$ and every $s\in\R_+$.
\item
The ring $\underline B$ satisfies condition
$\mathrm{(a)}^\mathrm{un}_\bff$ of \eqref{subsec_badabum}.
\end{enumerate}
\end{proposition}
\begin{proof}(i): Denote by $\Gamma_0\subset\Gamma$ the
submonoid generated by
$(\gamma_1/p^n,\dots,\gamma_r/p^n~|~n\in\N)$, and let
$\underline B{}_0:=\Gamma_0\times_\Gamma\underline B$.
Recall that $R_{r,0}:=\Z[P_r]$, where $P_r:=\N[1/p]^{\oplus r}$,
and let $\phi:P_r\to\Gamma_0$ be the unique morphism of
monoids such that $\phi(e_i):=\gamma_i$ for $i=1,\dots,r$,
where $(e_1,\dots,e_r)$ is the standard minimal system of
generators of $\N^{\oplus r}$. With this notation, $R_{r,0}$
is a $P_r$-graded $\Z$-algebra, and we set
$\underline S:=(R_{r,0})_{/\Gamma_0}$ (see definition
\ref{def_Gamma-graded-algs}(v)). Then clearly
$u$ induces a morphism of $\Gamma_0$-graded $\Z$-algebras
$\underline S\to\underline B{}_0$, and $\underline B$
can be regarded naturally as a $\Gamma$-graded
$\underline S{}_{/\Gamma}$-module. We remark, quite generally :

\begin{claim}\label{cl_slow-start}
Let $(\Delta,0,+)$ be a monoid, $\Delta_0\subset\Delta$ a
submonoid such that the translation maps $\Delta\to\Delta$ :
$\delta\mapsto\delta+\delta_0$ are injective for every
$\delta_0\in\Delta_0$. Let also $\underline S=(S,\gr_\bullet S)$
be a $\Delta_0$-graded $\Z$-algebra, $N$ a $\Delta_0$-graded
$\underline S$-module, and $M$ a $\Delta$-graded
$\underline S{}_{/\Delta}$-module. Then we have :
\begin{enumerate}
\item
The product $M':=\prod_{\delta\in\Delta}\gr_\delta M$
carries a natural $S$-module structure such that
the natural map $M\to M'$ is $S$-linear.
\item
The induced map $\tau_i:\Tor_i^S(N,M)\to\Tor_i^S(N,M')$
is injective, for every $i\in\N$.
\end{enumerate}
\end{claim}
\begin{pfclaim}(i): We define a $\Z$-bilinear map
$$
\mu_{\delta_0}:\gr_{\delta_0}S\times M'\to M'
\qquad
\text{for every $\delta_0\in\Delta_0$}
$$
by the rule : $(s,(m_\delta~|~\delta\in\Delta))\mapsto
(m'_\delta~|~\delta\in\Delta)$, where
$m'_\delta:=s\cdot m_\gamma$ if there exists a (necessarily
unique) $\gamma\in\Delta$ with $\gamma+\delta_0=\delta$,
and otherwise $m'_\delta:=0$. The sum of the maps $\mu_{\delta_0}$
yields a well defined $\Z$-bilinear map
$$
\mu:S\times M'\to M'
$$
and it is easily seen that $\mu$ defines an $S$-module
structure on $M'$ with the sought property.

(ii): In view of remark \ref{rem_graded-resol}(i), we may
find a resolution $L_\bullet\to N$ by free $S$-modules such
that $L_i$ is a $\Delta_0$-graded $\underline S$-module,
and the differential $d_{i+1}:L_{i+1}\to L_i$ is a morphism
of $\Delta_0$-graded $\underline S$-modules, for every
$i\in\N$. Then, for every $i\in\N$ the tensor product
$P_i:=L_i\otimes_SM$ carries a natural structure of
$\Delta$-graded $\underline S{}_{/\Delta}$-module : namely
$$
\gr_\delta P_i:=\Img\Bigl(\bigoplus_{(\delta_0,\gamma)\in D(\delta)}
\gr_{\delta_0}L_i\otimes_\Z\gr_\gamma M\to P_i\Bigr)
\qquad
\text{for every $\delta\in\Delta$}
$$
where $D(\delta)$ is the set of
all pairs $(\delta_0,\gamma)\in\Delta_0\times\Delta$ such
that $\gamma+\delta_0=\delta$. With this definition, it is
then easily seen that the differential
$d^P_{i+1}:=d_{i+1}\otimes_S\one_M$ of $P_\bullet$ is a
morphism of $\Delta$-graded
$\underline S{}_{/\Delta}$-modules, whence an induced
$\Delta$-graded $\underline S{}_{/\Delta}$-module
structure on $\Tor^S_i(N,M)$, for every $i\in\N$.
Moreover, we get a $\Z$-bilinear map
$$
\psi:L_i\times M'\to P'_i:=\prod_{\delta\in\Delta}\gr_\delta P_i
\qquad
\text{for every $i\in\N$}
$$
namely, the unique one whose restriction
$\gr_{\delta_0}L_i\times M'\to P'_i$ for every
$\delta_0\in\Delta_0$ is given by the rule :
$(l,(m_\delta~|~\delta\in\Delta))\mapsto
(l\otimes m'_\delta~|~\delta\in\Delta)$, where
$m'_\delta:=m_\gamma$ if $\gamma+\delta_0=\delta$, and
$m'_\delta:=0$ if there is no such $\gamma$. Then it
is easily seen that $\psi$ is even $S$-bilinear,
so it induces an $S$-linear map
$$
L_i\otimes_SM'\to P'_i
\qquad
\text{for every $i\in\N$}.
$$
Furthermore, the resulting diagram of $S$-linear maps
$$
\xymatrix{ L_{i+1}\otimes_SM' \ar[r]
\ar[d]_{d_{i+1}\otimes_S\one_{M^\wedge}} &
P'_{i+1} \ar[d]^{\prod_{\delta\in\Delta}\gr_\delta d^P_{i+1}} \\
L_i\otimes_SM' \ar[r] & P'_i
}$$
commutes for every $i\in\N$, and we deduce an $S$-linear map
$$
\Tor_i^S(N,M')\to\prod_{\delta\in\Delta}\gr_\delta\Tor_i^S(N,M)
\qquad
\text{for every $i\in\N$}
$$
whose composition with $\tau_i$ is the natural injection.
The claim follows.
\end{pfclaim}

Now, set $B':=\prod_{\gamma\in\Gamma}\gr_\gamma B$, endow
$B'$ with the linear topology defined as in remark
\ref{rem_graded-top-algs}(ii), and recall that the
inclusion map $j:B\to B'$ factors through the
inclusion map $B\to A$ and a continuous map $j_A:A\to B'$
(remark \ref{rem_graded-top-algs}(iii)). Endow $B'$
with the $S$-module structure provided by claim
\ref{cl_slow-start}(i); a direct inspection of the proof
of {\em loc.cit.} shows that scalar multiplication by any
homogeneous element $s$ of $S$ defines a continuous
endomorphism of $B'$, and on the other hand we have
$j(s\cdot b)=s\cdot j(b)$ for every such $s$ and every
$b\in B$; since $j$ is also a continuous map, we deduce
that the same identity holds more generally for every
$b\in A$, {\em i.e.} $j_A$ is an $S$-linear map. Lastly,
notice that the ideal $I^{\lceil s\rceil}_{r,0}$ of $\underline S$
is $\Gamma_0$-graded, so $R_{r,0}/I^{\lceil s\rceil}_{r,0}$ is a
$\Gamma_0$-graded $\underline S$-module for every $s\in\R_+$.
By claim \ref{cl_slow-start}, it follows that the composition
of the two natural maps
$$
\Tor_i^{R_{r,0}}(R_{r,0}/I^{\lceil s\rceil}_{r,0},\underline B)\to
\Tor_i^{R_{r,0}}(R_{r,0}/I^{\lceil s\rceil}_{r,0},A)\to
\Tor_i^{R_{r,0}}(R_{r,0}/I^{\lceil s\rceil}_{r,0},B')
$$
is injective for every $i\in\N$; on the other hand, the
middle term vanishes for $i>0$, due to proposition
\ref{prop_back-to-toids}, whence the contention.

(ii) follows from (i) and lemma \ref{lem_generalissimo}.
\end{proof}

\sset\subsubsection{}\label{subsec_real-blowup}
We return now to the situation of
\eqref{subsec_general-case}, and we suppose we are given
a family $(\cK_\gamma~|~\gamma\in\Gamma)$ consisting of
topologically closed $1$-taut bounded $\bar\bE$-submodules
of $\bE_U$, indexed by a $p$-perfect monoid $(\Gamma,+,0)$,
and such that :
\begin{enumerate}
\alphaenu
\item
$1\in\cK_0$
\item
$\cK_{p\gamma}=\cK_\gamma^p$ for every $\gamma\in\Gamma$
\item
$\cK_\gamma\cdot\cK_\mu\subset\cK_{\gamma+\mu}$
for every $\gamma,\mu\in\Gamma$
\end{enumerate}
(see \eqref{subsec_quasi-taut} and \eqref{subsec_powers-and-Witt},
and notice that a subset of $\bE_U$ is bounded if and only
if it is contained in $\Fil^s\bE_U$, for some $s\in\Z[1/p]$).
We set
$$
\cE:=\bigoplus_{\gamma\in\Gamma}\cK_\gamma
\qquad
\cA:=\bigoplus_{\gamma\in\Gamma}\{\cK_\gamma\}.
$$
It is easily seen that the multiplication of $\bE_U$ (resp.
of $A_U$) induces a $\Gamma$-graded $\bar\bE$-algebra
structure (resp. $\bar A$-algebra structure) on $\cE$
(resp. on $\cA$). Choose ideals of definition $J$ for
$A$ and $\cJ$ for $\bE$, and endow $\cE$ (resp. $\cA$)
with its $\cJ$-adic (resp. $J$-adic) topology. Let also
$\cE^\wedge$ (resp. $\cA^\wedge$) be the separated completion
of $\cE$ (resp. of $\cA$).

\begin{lemma} The pairs
$$
(\cE{}^\wedge,\cE)
\qquad\text{and}\qquad
(\cA{}^\wedge,\cA)
$$
are topological rings with $\Gamma$-graded structures.
\end{lemma}
\begin{proof} Indeed, the assertion comes down to
checking that $\cK_\gamma$ (resp. $\{\cK_\gamma\}$)
is a closed subset of $\cE{}^\wedge$ (resp. of
$\cA{}^\wedge$) for every $\gamma\in\Gamma$, and notice
that the topological closure of $\cK_\gamma$ in $\cE^\wedge$
(resp. of $\{\cK_\gamma\}$ in $\cA^\wedge$) is the $\cJ$-adic
(resp. $J$-adic) completion of $\cK_\gamma$ (resp. of
$\{\cK_\gamma\}$). Now, by assumption $\cK_\gamma$ is bounded
in $\bE_U$, so $\{\cK_\gamma\}$ is bounded in $A_U$, for every
$\gamma\in\Gamma$; it follows that the topology $\cT_\gamma$
on $\cK_\gamma$ (resp. on $\{\cK_\gamma\}$) induced by the
inclusion into $\bE_U$ (resp. into $A_U$) is coarser than
the $\cJ$-adic (resp. $J$-adic) topology; on the other
hand, $\cT_\delta$ is complete and separated, since the
same holds for the topology of $\bE_U$ (resp. of $A_U$).
Then lemma \ref{lem_fontaine} implies that $\cK_\gamma$
(resp. $\{\cK_\gamma\}$) is complete and separated for its
$\cJ$-adic (resp. $J$-adic) topology, whence the contention.
\end{proof}

\begin{proposition}\label{prop_new-graded-perf}
With the notation of \eqref{subsec_real-blowup}, the
following holds :
\begin{enumerate}
\item
The rings $\cE{}^\wedge$, and $\cA{}^\wedge$ are perfectoid.
\item
There is a natural isomorphism of topological rings with
$\Gamma$-graded structures
$$
\bE(\cA{}^\wedge,\cA)\isom(\cE{}^\wedge,\cE).
$$
\end{enumerate}
\end{proposition}
\begin{proof} Taking into account remark
\ref{rem_topology-of-E}(v) and example
\ref{ex_discrete-Witt}(ii), it is easily seen that
$\cE{}^\wedge$ is a perfect topological $\bE$-algebra, hence
it is perfectoid (see example \ref{ex_perfectoid}(i)). Set
$$
(\cW^\wedge,\cW):=W(\cE{}^\wedge,\cE)
\qquad
(\bar\cW{}^\wedge,\bar\cW):=(\cW^\wedge,\cW)\otimes_{W(\bE)}A.
$$
To prove both assertions, it will therefore suffice to
exhibit isomorphisms
\set\begin{equation}\label{eq_too-late-for-yoga}
(\bar\cW{}^\wedge,\bar\cW)\isom(\cA{}^\wedge,\cA)
\end{equation}
of topological $A$-algebras with $\Gamma$-graded structures.
We shall first exhibit a continuous map
$$
u_\cA:(\cW{}^\wedge,\cW)\to(\cA{}^\wedge,\cA).
$$
To this aim, notice first that condition (b) of
\eqref{subsec_real-blowup} implies that
$$
\gr_\gamma\cW=W(\cK_\gamma)
\qquad
\text{for every $\gamma\in\Gamma$}
$$
(notation of \eqref{subsec_powers-and-Witt}) and the
multiplication law of $\cW$ is the unique bilinear
map $\cW\times\cW\to\cW$ whose restriction
$$
W(\cK_\gamma)\times W(\cK_\mu)\to
W(\cK_{\delta+\gamma})
$$
is induced by the multiplication law of $W(\bE_U)$, for every
$\gamma,\mu\in\Gamma$. Also, since $\cK_\gamma$ is closed and
$1$-taut, lemma \ref{lem_like-th-taut-two}(iii) says that $u_U$
restricts to a well defined surjective map
$$
\gr_\gamma u_\cA:W(\cK_\gamma)\to\gr_\gamma\cA
\qquad
\text{for every $\gamma\in\Gamma$}
$$
whose kernel is $\underline\alpha W(\cK_\gamma)$, for any
distinguished element $\underline\alpha$ of $\Ker\,u_A$.
Summing up, we conclude that the direct sum of the
system of maps $(\gr_\gamma u_\cA~|~\gamma\in\Gamma)$
is a well defined surjective ring homomorphism
$\gr_\bullet u_\cA$ that makes commute the diagram
\set\begin{equation}\label{eq_simpler-than-before}
{\diagram
W(\bE) \ar[r]^-{u_A} \ar[d] & A \ar[d] \\
\cW \ar[r]^{\gr_\bullet u_\cA} & \cA
\enddiagram}
\end{equation}
and with kernel equal to $\underline\alpha\cW$. Moreover,
it follows easily from propositions \ref{prop_morel}(ii)
and \ref{prop_general-graded-top}(i.b) that the left
vertical map of \eqref{eq_simpler-than-before} is an adic
ring homomorphism, and the same obviously holds for the
right vertical map. Since $u_A$ is an open map, we deduce
that $\gr_\bullet u_\cA$ is continuous and open as well;
especially, it extends to a well defined continuous
surjective ring homomorphism $u_\cA$ as sought, whose kernel
is the topological closure of $\underline\alpha\cW$ in
$\cW^\wedge$ (proposition \ref{prop_replaces-Mat-Th.8.1}).
However, recall that $\underline\alpha\cW^\wedge$ is a closed
ideal of $\cW^\wedge$ (proposition \ref{prop_bar-W-perfect}(ii)),
and obviously $\underline\alpha\cW$ is dense in it, hence
the kernel of $u_\cA$ is $\underline\alpha\cW^\wedge$, as
required. Lastly, since $\gr_\bullet u_\cA$ is an open map,
the same holds for $u_\cA$, and the proof is concluded.
\end{proof}

\begin{example}\label{ex_E-of-beta}
In the situation of \eqref{subsec_general-case}, choose
ideals of definition $J$ for $A$, and $\cJ$ for $\bE$.
Let also $\beta_\bullet:=(\beta_\lambda~|~\lambda\in\Lambda)$
be any bounded system of elements of $\bE_U$. We consider
the unique morphism of monoids
$$
\phi_\bE:P:=\N[1/p]^{(\Lambda)}\to\bE_U
\qquad
e_\lambda\mapsto\beta_\lambda
\qquad
\text{for every $\lambda\in\Lambda$}
$$
where $e_\bullet:=(e_\lambda~|~\lambda\in\Lambda)$ is the
standard minimal set of generators of $\N^{(\Lambda)}$.
Let also $I\subset P$ be the ideal generated by $e_\bullet$.
The composition $\phi_A:\phi^\flat_U\circ\phi_\bE:P\to A$
is also a morphism of monoids, and we let
$I^{\La\delta\Ra}\bar\bE$ (resp. $I^{\La\delta\Ra}\bar A$) be the
$\bar\bE$-submodule of $\bE_U$ generated by
$\phi_\bE(I^{\La\delta\Ra})$ (resp. the $\bar A$-submodule
of $A_U$ generated by $\phi_A(I^{\La\delta\Ra})$), for every
$\delta\in\Delta:=\N[1/p]$ (notation of
\eqref{subsec_mon-fract-powers}).
We denote by $\gr_\delta\cE$ (resp. $\gr_\delta\cA$)
the topological closure of $I^{\La\delta\Ra}\bar\bE$ in
$\bE_U$ (resp. of $I^{\La\delta\Ra}\bar A$ in $A_U$) for
every $\delta\in\Delta$ (notation of remark
\ref{rem_Witt-are-f-adic}(i)). We set
$$
\cE(\beta_\bullet):=\bigoplus_{\delta\in\Delta}\gr_\delta\cE
\qquad
\cA(\beta_\bullet):=\bigoplus_{\delta\in\Delta}\gr_\delta\cA.
$$
Hence, $\cE(\beta_\bullet)$ is a $\Delta$-graded $\bE$-algebra,
and $\cA(\beta_\bullet)$ is a $\Delta$-graded $A$-algebras. We
endow $\cE(\beta_\bullet)$ with the $\cJ$-adic topology and
$\cA(\beta_\bullet)$ with the $J$-adic topology, and let
$\cE(\beta_\bullet)^\wedge$ and $\cA(\beta_\bullet)^\wedge$ be
the respective separated completion. It is easily seen that
$\gr_\delta\cE$ is $1$-taut and bounded in $\bE_U$, and moreover
$\gr_{p\delta}\cE=(\gr_\delta\cE)^p$ for every $\delta\in\Delta$;
furthermore, $\gr_\delta\cA=\{\gr_\delta\cE\}$ for every such
$\delta$. Then proposition \ref{prop_new-graded-perf} shows
that $(\cE(\beta_\bullet)^\wedge,\cE(\beta_\bullet))$, and
$(\cA(\beta_\bullet)^\wedge,\cA(\beta_\bullet))$ are perfectoid
rings with $\Delta$-graded structures, and we have a natural
isomorphism of topological rings with $\Delta$-graded structures :
$$
\bE(\cA(\beta_\bullet)^\wedge,\cA(\beta_\bullet))\isom
(\cE(\beta_\bullet)^\wedge,\cE(\beta_\bullet)).
$$
\end{example}

\begin{example}\label{ex_filtered-type}
In the situation of \eqref{eq_back-on-track}, let
$\Gamma:=\Z[1/p]$; the family
$(\Fil^\gamma\bE_U~|~\gamma\in\Gamma)$ fulfills conditions
(a)--(c) of \eqref{subsec_real-blowup}, and taking into
account corollary \ref{cor_boody-boody} and proposition
\ref{prop_new-graded-perf}, we get perfectoid rings with
$\Gamma$-graded structures $(\cE_U^\wedge,\cE_U)$ and
$(\cA_U^\wedge,\cA_U)$, by setting
$$
\cE_U:=\bigoplus_{\gamma\in\Gamma}\Fil^\gamma\bE_U
\qquad
\cA_U:=\bigoplus_{\gamma\in\Gamma}\Fil^\gamma A_U
$$
as well as a natural isomorphism of topological rings with
$\Gamma$-graded structures
\set\begin{equation}\label{eq_vladimir}
\bE(\cA_U^\wedge,\cA_U)\isom(\cE_U^\wedge,\cE_U).
\end{equation}
Moreover, by inspecting the proof of proposition
\ref{prop_new-graded-perf} it is easily seen that --
under the identification \eqref{eq_vladimir} -- the
restriction $\gr_\gamma\cE_U\to\gr_\gamma\cA_U$ of the
map $\bar u_{\cA^\wedge_U}:\cE_U^\wedge\to\cA^\wedge_U$ agrees
with the restriction $\Fil^\gamma\bE_U\to\Fil^\gamma A_U$
of the map $\phi^\flat_U:\bE_U\to A_U$, for every
$\gamma\in\Gamma$.
\end{example}

\sset\subsubsection{}\label{subsec_approx-with-Teich}
Let $\Gamma$ be a $p$-perfect monoid, $(A^\wedge,A)$ a
perfectoid ring with $\Gamma$-graded structure, set
$(\bE^\wedge,\bE):=\bE(A^\wedge,A)$ and let
$\underline\alpha\in W(\gr_0\bE)$ be a distinguished
element of $\Ker\,u_A$ (see proposition
\ref{prop_graded-perfectoid}(iii)). We denote by
$\bar u_A:\bE\to A$ the restriction of the map
$u_{A^\wedge}:\bE^\wedge\to A^\wedge$. For every integer
$n>0$ let
$$
\lambda_n:=\sum_{i=0}^{n-1}p^{-i}
\qquad
S_n:=\{(i,j)\in\N[1/p]^{\oplus 2}~|~\text{$\lambda_n\leq i<n$
\ \ \text{and}\ \ $1>j>1-i/n$}\}
$$
and for every $\gamma\in\Gamma$ and every $e\in\gr_\gamma\bE$,
define the $\gr_0A$-module
$$
\cI(e,n,\gamma):=
\sum_{(i,j)\in S_n}\gr_{(1-j)\gamma}A\cdot\bar u_A(\alpha_0^i\cdot e^j).
$$

\begin{proposition}\label{prop_Scholze-approx-lemma}
With the notation of \eqref{subsec_approx-with-Teich}, let
$\gamma\in\Gamma$ and $a\in\gr_\gamma A$ be any elements. Then,
for every integer $n>0$ there exists $e\in\gr_\gamma\bE$ such
that
$$
a-\bar u_A(e)\in\cI(e,n,\gamma)+p^n\cdot\gr_\gamma A.
$$
\end{proposition}
\begin{proof} Notice that the assertion is independent of
the topology $\cT_{A^\wedge}$ of $A^\wedge$, so we may assume
that $\cT_{A^\wedge}$ is the $p$-adic topology of $A^\wedge$
(proposition \ref{prop_change-topol}(i)). We argue by
induction on $n\in\N$. For $n=1$, we have $S_1=\emptyset$,
and the assertion means that there exists $e\in\gr_\gamma\bE$
such that $a-\bar u_A(e)\in p\cdot\gr_\gamma A$. Hence, let
$\pi:A\to A/pA$ be the projection; it suffices to remark :

\begin{claim}
The map $\pi\circ\bar u_A:\bE\to A/pA$ is surjective.
\end{claim}
\begin{pfclaim} The assertion follows easily from
claim \ref{cl_graded-frob}(i).
\end{pfclaim}

Thus, suppose that the proposition has already been proven
for some $n\geq 1$, for every $\gamma\in\Gamma$ and every
$a\in\gr_\gamma A$. We notice :

\begin{claim}\label{cl_Z-linear-comb}
(i)\ \
For every $\delta\in\Gamma$, every $b\in\gr_\delta A$ and
every $k\in\N$ there exist a sequence
$(\beta_0,\dots,\beta_k)$ of elements of $\gr_\delta\bE$,
and an element $d\in\gr_\delta A$ such that
\set\begin{equation}\label{eq_easy-redux}
b=\sum_{i=0}^kp^i\cdot\bar u_A(\beta_i)+p^{k+1}d.
\end{equation}

(ii)\ \
There exist elements $c\in\gr_\gamma A$, $e_n\in\gr_\gamma\bE$,
a finite set $\Lambda$, a system of positive integers
$(k_\lambda~|~\lambda\in\Lambda)$, a mapping
$\Lambda\to S_n$ : $\lambda\mapsto(i_\lambda,j_\lambda)$,
a system $(f_\lambda~|~\lambda\in\Lambda)$ with
$f_\lambda\in\gr_{(1-j_\lambda)\gamma}\bE$ for every
$\lambda\in\Lambda$, such that
\set\begin{equation}\label{eq_Z-linear-comb}
a=\bar u_A(e_n)+\sum_{\lambda\in\Lambda}k_\lambda\cdot
\bar u_A(f_\lambda\cdot\alpha_0^{i_\lambda}\cdot e_n^{j_\lambda})
+p^nc.
\end{equation}

(iii)\ \
There exists $f\in\gr_0\bE$ such that
$p^n-\bar u_A(\alpha_0\cdot f)^n\in p^{n+1}\gr_0A$
for every $n\in\N$.
\end{claim}
\begin{pfclaim}(i): We argue by induction on $k$,
and notice that the assertion for $k=0$ follows
from the case $n=1$ of the proposition, which is
already known. Thus, suppose that $k\in\N$, and
that we have already found a sequence
$(\beta_0,\dots,\beta_k)$ such that \eqref{eq_easy-redux}
holds. Then we may also find $\beta_{k+1}\in\gr_\delta\bE$
and $d'\in\gr_\delta A$ such that $d=\bar u_A(\beta_{k+1})+pd'$,
and replacing this expression for $d$ in \eqref{eq_easy-redux}
we get the sought identity for $k+1$.

(ii): By inductive assumption, there exist a finite
subset $T\subset S_n$ and for every $(i,j)\in T$ an
element $b_{(i,j)}\in\gr_{(1-j)\gamma}A$ such that
$$
a=\bar u_A(e_n)+\sum_{(i,j)\in T}b_{(i,j)}\cdot
\bar u_A(\alpha_0^i\cdot e_n^j)+p^nd
\qquad
\text{for some $e_n\in\gr_\gamma\bE$ and $d\in\gr_\gamma A$}.
$$
By (i), we may then find, for every $(i,j)\in T$, a
system $(f_{(i,j),0},\dots,f_{(i,j),n-1})$ of elements of
$\gr_{(1-j)\gamma}\bE$ and an element $c_{(i,j)}\in\gr_{(1-j)\gamma}A$
such that
$b_{(i,j)}=\sum_{i=0}^{n-1}p^i\cdot\bar u_A(f_{(i,j),i})+p^nc_{(i,j)}$.
Hence, the sought identity holds with
$$
c:=d+\!\!\sum_{(i,j)\in T}c_{(i,j)}\cdot\bar u_A(\alpha_0^i\cdot e_n^j)
\qquad
\Lambda:=T\times\{0,\dots,n-1\}
\qquad
k_{(t,s)}:=p^s
\quad
\text{for every $(t,s)\in\Lambda$}
$$
and with the mapping $\Lambda\to S_n$ : $(t,k)\mapsto t$.

(iii): By lemma \ref{lem_was-third-cond}(iii), there exists
$x\in A^\wedge$ such that $p=\bar u_A(\alpha_0)\cdot x$. Let
$\pi^\wedge_0:A^\wedge\to\gr_0A$ be the canonical $0$-projection
(see remark \ref{rem_graded-top-algs}(iii)); it follows
that $p=\bar u_A(\alpha_0)\cdot\pi^\wedge_0(x)$, so we may
assume that $x\in\gr_0A$. By (i), we have
$y_0:=x-\bar u_A(f)\in p\cdot\gr_0A$ for some $f\in\bE$,
and we check, by induction on $n$, that such $f$ will do.
Indeed, for $n=0$ the assertion is trivial, and for $n=1$
we have
$$
p=\bar u_A(\alpha_0)\cdot(y_0+\bar u_A(f))=
\bar u_A(\alpha_0\cdot f)+y_0\cdot\bar u_A(\alpha_0)
$$
and $y_1:=y_0\cdot\bar u_A(\alpha_0)=
py_0\cdot\pi_0^\wedge(x^{-1})\in p^2\cdot\gr_0A$, again
by lemma \ref{lem_was-third-cond}(iii). Suppose now
that $n\geq 1$, and that we have already shown that
$y_n:=p^n-\bar u_A(\alpha_0\cdot f)^n\in p^{n+1}\gr_0A$.
We compute
$$
\begin{aligned}
py_n=\,&
p^{n+1}-p\cdot\bar u_A(\alpha_0\cdot f)^n \\
=\,&
p^{n+1}-(y_1+\bar u_A(\alpha_0\cdot f))
\cdot\bar u_A(\alpha_0\cdot f)^n \\
=\,&
p^{n+1}-\bar u_A(\alpha_0\cdot f)^{n+1}-
y_1\cdot\bar u_A(\alpha_0\cdot f)^n
\end{aligned}
$$
Since both $py_n$ and $y_1\cdot\bar u_A(\alpha_0\cdot f)^n=
p^ny_1\cdot\bar u_A(f)^n\cdot\pi^\wedge_0(x^{-n})$ lie in
$p^{n+2}\gr_0A$, the assertion follows for $n+1$.
\end{pfclaim}

Pick any identity \eqref{eq_Z-linear-comb} as provided
by claim \ref{cl_Z-linear-comb}(ii), and any $f\in\gr_0\bE$
such that claim \ref{cl_Z-linear-comb}(iii) holds; by claim
\ref{cl_Z-linear-comb}(i) we may also find $g\in\gr_\gamma\bE$
such that $c-\bar u_A(g)\in p\cdot\gr_\gamma A$, and notice
that, after replacing the set $\Lambda$ by
$\bigcup_{\lambda\in\Lambda}\{\lambda\}\times\{1,\dots,k_\lambda\}$,
we may also assume that $k_\lambda=1$ for every
$\lambda\in\Lambda$, so that
$$
a=\bar u_A(e_n)+\sum_{\lambda\in\Lambda}
\bar u_A(f_\lambda\cdot\alpha_0^{i_\lambda}\cdot e_n^{j_\lambda})
+\bar u_A(g\cdot\alpha_0^n\cdot f^n)
+p^{n+1}d
\qquad
\text{for some $d\in\gr_\gamma A$}.
$$
We let $S'_n:=\{(i,j)\in\N[1/p]^{\oplus 2}~|~
\text{$\lambda_{n+1}\leq i<n+1$,\ \ $nj+i\geq 1+n$
\ \ and\ \ $j<1$}\}$ and we set
$$ 
\delta:=\sum_{\lambda\in\Lambda}
f_\lambda\cdot\alpha_0^{i_\lambda}\cdot e_n^{j_\lambda}
\qquad\text{and}\qquad
e_{n+1}:=e_n+\delta+g\cdot\alpha_0^n\cdot f^n.
$$

\begin{claim}\label{cl_new-expression}
There exists a finite subset $\Lambda'$ with a mapping
$\Lambda'\to S'_n$ : $\lambda\mapsto(i'_\lambda,j'_\lambda)$,
and for every $\lambda\in\Lambda'$ an element
$z_\lambda\in\gr_{(1-j'_\lambda)\gamma}A$ such that
$$
a=\bar u_A(e_{n+1})+\sum_{\lambda\in\Lambda'}
z_\lambda\cdot\bar u_A(\alpha_0^{i'_\lambda}\cdot e^{j'_\lambda}_n)
+p^{n+1}d'
\qquad
\text{for some $d'\in\gr_\gamma A$}.
$$
\end{claim}
\begin{pfclaim} For every $r=1,\dots,n$, let $\Lambda'_r$
be the subset of all sequences
$\underline\sigma:=
(\sigma_0,\sigma_1,\sigma_\lambda~|~\lambda\in\Lambda)$
of elements of $p^{-r}\N\setminus p^{1-r}\N$ such that
$\sigma_0+\sigma_1+\sum_{\lambda\in\Lambda}\sigma_\lambda=1$,
and set $\Lambda':=\bigcup_{r=1}^n\Lambda'_r$. Let also
$x:=a-\bar u_A(e_{n+1})$; according to proposition
\ref{prop_combinatorial}, there exist a mapping
$\Lambda'\to\Z_p$ : $\underline\sigma\mapsto z'_{\underline\sigma}$
and an element $d'\in\gr_\gamma A$ such that
$$
\begin{aligned}
x=&\, 
\sum_{r=1}^np^r\cdot\sum_{\underline\sigma\in\Lambda'_r}z'_{\underline\sigma}\cdot
\bar u_A\Bigl(e_n^{\sigma_0}\cdot g^{\sigma_1}\cdot f^{n\sigma_1}
\cdot\alpha_0^{n\sigma_1}\cdot\prod_{\lambda\in\Lambda}f_\lambda^{\sigma_\lambda}
\cdot\alpha_0^{\sigma_\lambda i_\lambda}\cdot e_n^{\sigma_\lambda j_\lambda}\Bigr)
+p^{n+1}d' \\
=&\,
\sum_{r=1}^np^r\cdot\sum_{\underline\sigma\in\Lambda'_r}z'_{\underline\sigma}\cdot
\bar u_A(f'_{\underline\sigma}\cdot\alpha_0^{i''_{\underline\sigma}}\cdot
e_n^{j'_{\underline\sigma}})
+p^{n+1}d'
\end{aligned}
$$
where
$$
f'_{\underline\sigma}:=
g^{\sigma_1}\cdot f^{n\sigma_1}\cdot\prod_{\lambda\in\Lambda}
f_\lambda^{\sigma_\lambda}
\qquad
i''_{\underline\sigma}:=
n\sigma_1+\sum_{\lambda\in\Lambda}\sigma_\lambda i_\lambda
\qquad
j'_{\underline\sigma}:=
\sigma_0+\sum_{\lambda\in\Lambda}\sigma_\lambda j_\lambda
\qquad
\text{for every $\underline\sigma\in\Lambda'$}.
$$
Let us remark that
\set\begin{equation}\label{eq_greater-than-one}
j'_{\underline\sigma}+\frac{i''_{\underline\sigma}}{n}=
\sigma_0+\sigma_1+
\sum_{\lambda\in\Lambda}\sigma_\lambda\cdot
\Bigl(j_\lambda+\frac{i_\lambda}{n}\Bigr)\geq
\sigma_0+\sigma_1+\sum_{\lambda\in\Lambda}\sigma_\lambda=1
\quad
\text{for every $\underline\sigma\in\Lambda'$}.
\end{equation}
Moreover, we notice that
\set\begin{equation}\label{eq_second-bound}
1>j'_{\underline\sigma}
\qquad\text{and}\qquad
n>i''_{\underline\sigma}\geq\frac{\lambda_n}{p^r}
\qquad
\text{for every $r=1,\dots,n$ and every
$\underline\sigma\in\Lambda'_r$}.
\end{equation}
Indeed, the first inequality is immediate, since $\sigma_0<1$;
the second inequality is clear if $\sigma_\lambda\neq 0$ for
at least one index $\lambda\in\Lambda$, since in this case
$\sigma_\lambda\geq p^{-r}$ and
$i''_{\underline\sigma}<n\sigma_1+n\sum_{\lambda\in\Lambda}\sigma_\lambda<n$.
If $\sigma_\lambda=0$ for every $\lambda\in\Lambda$, we must have
$1>\sigma_1>0$, since otherwise we would get
$\sigma_0,\sigma_1\in\N$, which is absurd; hence, in this case
$1>\sigma_1\geq p^{-r}$, so
$n>i''_{\underline\sigma}\geq n/p^r\geq\lambda_n/p^r$.

Next, for every $r=1,\dots,n$ and every
$\underline\sigma\in\Lambda'_r$, we may find
$z_{\underline\sigma}\in\gr_{(1-j'_{\underline\sigma})\gamma}A$ such that
$$
p^r\cdot z'_{\underline\sigma}\cdot\bar u_A(f'_{\underline\sigma})=
z_{\underline\sigma}\cdot\bar u_A(\alpha_0^r)
$$
(cp. the proof of claim \ref{cl_Z-linear-comb}(iii)).
Hence, if we set
$$
i'_{\underline\sigma}:=i''_{\underline\sigma}+r
\qquad
\text{for every $r=1,\dots,n$ and every
$\underline\sigma\in\Lambda'_r$}
$$
we get an expression for $a$ of the sought type,
and it remains only to check that
$(i'_{\underline\sigma},j'_{\underline\sigma})\in S'_n$ for
every $\underline\sigma\in\Lambda'$. However, from
\eqref{eq_second-bound} we see that
$i'_{\underline\sigma}\geq 1+\lambda_n/p=\lambda_{n+1}$
if $\underline\sigma\in\Lambda'_1$, and
$i'_{\underline\sigma}\geq 2+\lambda_n/p^2>\lambda_{n+1}$
if $\underline\sigma\in\Lambda'_r$ for some $r>1$.
Lastly, from \eqref{eq_greater-than-one} we get
$i''_{\underline\sigma}+j'_{\underline\sigma}\cdot n\geq n$
for every $\underline\sigma\in\Lambda'$, whence
$i'_{\underline\sigma}+j'_{\underline\sigma}\cdot n\geq n+1$,
and the proof is complete.
\end{pfclaim}

Let $I$ (resp. $J$) be the ideal of $\bE$ generated by the
system $(\alpha^n_0,e_n)$ (resp. by $(\alpha_0^n,e_n+\delta)$) and
set $q:=\min(j+i/n~|~(i,j)\in\Lambda)$; then
$\delta\in\mathrm{i.c.}(I,\bE,q)$, and the radical of $I$ equals
the radical of $J$. Since $q>1$, lemma \ref{lem_return-to-ic}(ii)
shows therefore that
$$
I^{\La 1\Ra}\bE=J^{\La 1\Ra}\bE
$$
and consequently $e_n^t\in J^{\La t\Ra}\bE$ for every $t\in\N[1/p]$.
Notice that $(\alpha_0^n,e_{n+1})$ is another system of
generators for $J$. Since the topology of $A$ is $p$-adic,
it follows that for every $t\in\N[1/p]$ there exist a finite
subset $T_t\subset\N[1/p]$ with $0\leq\sigma\leq t$ for every
$\sigma\in T_t$, and a system $(h'_\sigma~|~\sigma\in T_t)$
of elements of $A^\wedge$ such that
$$
\bar u_A(e^t_n)=
\sum_{\sigma\in T_t}h'_\sigma\cdot
\bar u_A(\alpha_0^{n\sigma}\cdot e_{n+1}^{t-\sigma})
$$
(corollary \ref{cor_taut-two}(ii)); after replacing
$h'_\sigma$ with its canonical $\sigma\gamma$-projection,
we may also assume that $h'_\sigma\in\gr_{\sigma\gamma}A$ for
every $\sigma\in T_t$. Now, pick an expression for $a$ as
in claim \ref{cl_new-expression}; to conclude the proof,
it suffices to show that
$\bar u_A(\alpha_0^{i'_\lambda}\cdot e_n^{j'_\lambda})\in
\cI(e_{n+1},n+1,\gamma)+p^{n+1}\cdot\gr_{j'_\lambda\gamma}A$
for every $\lambda\in\Lambda'$. To ease notation, set
$(s,t):=(i'_\lambda,j'_\lambda)$; we see that
$$
\bar u_A(\alpha_0^s\cdot e_n^t)=
\sum_{\sigma\in T_t} h'_\sigma\cdot
\bar u_A(\alpha_0^{s+n\sigma}\cdot e^{t-\sigma}_{n+1}).
$$
Now, since $(s,t)\in S'_n$, we have
$s+n\sigma\geq s\geq\lambda_{n+1}$ for every $\sigma\in T'_t$.
Moreover :
$$
s+n\sigma+(t-\sigma)\cdot(n+1)=s+nt+t-\sigma\geq n+1+t-\sigma.
$$
Recall also that $s+nt\geq n+1$; therefore, if $s+n\sigma<n+1$,
we must have $t-\sigma>0$, so the pair $(s+n\sigma,t-\sigma)$
lies in $S_{n+1}$; lastly, if $s+n\sigma\geq n+1$, the term
$\bar u_A(\alpha_0^{s+n\sigma}\cdot e^{t-\sigma}_{n+1})$
lies in $p^{n+1}\cdot\gr_{(t-\sigma)\gamma}A$, whence
$h'_\sigma\cdot\bar u_A(\alpha_0^{s+n\sigma}\cdot e^{t-\sigma}_{n+1})
\in p^{n+1}\cdot\gr_{t\gamma}A$ and the assertion follows.
\end{proof}

\begin{remark}\label{rem_Scholze-approx-lemma}
(i)\ \
One can show by a direct computation that the element
$e_n$ appearing in the proof of proposition
\ref{prop_Scholze-approx-lemma} is integral over the
ideal generated by $e_{n+1}$ and $\alpha_0^n$; as a
consequence, the proof of {\em loc. cit.} can be made
in principle completely effective, so that one can
extract from it an algorithm that exhibits an element
$e$ with the sought properties, for any given $a\in\gr_\gamma A$.

(ii)\ \
With a more careful tracking of the intervening estimates,
one can also prove a refinement of proposition
\ref{prop_Scholze-approx-lemma}, whereby the set of
allowed exponents $S_n$ is replaced by the slightly smaller
one, consisting of all pairs $(i,j)\in\N[1/p]^{\oplus 2}$
such that $n>i\geq\lambda_n$ and $j\geq f_n(i)$, where
$f_n:[\lambda_n,+\infty[\to\R$ is the piecewise linear
continuous function such that :
\begin{itemize}
\item
$f_n(x)=0$ for every $x\geq n$.
\item
$f_n$ has slope $-1/k$ on the segment
$[\lambda_{n-k}+k/p^{n-k},\lambda_{n-k-1}+(k+1)/p^{n-k-1}]$,
for every integer $k=1,\dots,n-1$.
\end{itemize}

(iii)\ \
In \cite{Scho}, Scholze proves theorem \ref{th_Scholze-tilt}
(for his perfectoid rings, that are a special case of our
perfectoid quasi-affinoid rings) by means of his
``approximation lemma'' (see \cite[Lemma 6.5]{Scho}). With
our proposition \ref{prop_Scholze-approx-lemma}, we can give
an alternative proof of theorem \ref{th_Scholze-tilt} that
generalizes to quasi-affinoid perfectoid rings, along the lines
of Scholze's original argument. The first step is the following
extension of proposition \ref{prop_Cont-u_A} to quasi-affinoid
perfectoid rings.
\end{remark}

\begin{lemma}\label{lem_tilt-Scholze}
In the situation of \eqref{subsec_general-case},
let $v$ be any continuous valuation on $A_U$. Then :
\begin{enumerate}
\item
The mapping $v\circ\phi^\flat_U$ is a continuous
valuation of\/ $\bE_U$.
\item
The rule : $v\mapsto v\circ\phi^\flat_U$ defines a continuous map
$$
\Cont(\phi^\flat_U):\Cont(A_U)\to\Cont(\bE_U).
$$
\end{enumerate}
\end{lemma}
\begin{proof}(i): Since both $v$ and $\phi^\flat_U$ are morphisms
of pointed monoids, the same holds for $v\circ\phi^\flat_U$. Thus,
it remains only to check that
$$
v(\phi^\flat_U(x+y))\leq
M(x,y):=\max(v\circ\phi^\flat_U(x),v\circ\phi^\flat_U(y))
\qquad
\text{for every $x,y\in\bE_U$}.
$$
However, from proposition \ref{prop_new-formula}(ii) we know
that $\phi^\flat_U(x+y)$ is the sum of
$\phi^\flat_U(x)+\phi^\flat_U(y)$ and a $p$-adically convergent
series, whose $n$-th term can be written as $p^n\cdot z_n$,
where $z_n$ is in turn a finite $\Z_p$-linear combination of
elements of the form $\phi^\flat_U(x^\sigma y^{1-\sigma})$, with
$\sigma$ ranging over the elements of
$\Sigma_n:=p^{-n}\N\cap[0,1]$. It follows that
$$
v(z_n)\leq\max
(\sigma\cdot v(\phi^\flat_U(x))+(1-\sigma)\cdot v(\phi^\flat_U(y))
~|~\sigma\in\Sigma_n)\leq M(x,y).
$$
On the other hand, $v(p)$ is cofinal in the value group of $v$,
since the latter is continuous (see \eqref{subsec_I-admissible}).
The assertion is an immediate consequence.

The proof of (ii) is similar to that of proposition
\ref{prop_Cont-u_A}(ii) : details left to the reader. 
\end{proof}

Next, we remark the following corollary of proposition
\ref{prop_Scholze-approx-lemma}, that generalizes and
improves slightly Scholze's ``approximation lemma'' in
\cite[Lemma 6.5]{Scho}. The same improvement appears in
\cite[Cor.3.6.7]{Ked-Liu}.

\begin{corollary}\label{cor_Scholze-approx-lemma}
Resume the notation of \eqref{eq_back-on-track}. For every
$\gamma\in\Z[1/p]$, every $a\in\Fil^\gamma A_U$, every $m\in\N$
and every $\beta\in\N[1/p]$ with $\beta<p/(p-1)$ there exists
$e\in\Fil^\gamma\bE_U$ such that
$$
v(a-\phi_U^\flat(e))\leq
v(p)^\beta\cdot\max(v(\phi^\flat_U(e)),v(p)^m)
\qquad
\text{for every $v\in\Cont(A_U)$}.
$$
\end{corollary}
\begin{proof} (Notice that the value group $\Gamma_v$ is
$p$-divisible, due to theorem \ref{th_Scholze-tilt}(iv),
but even discounting {\em loc.cit.}, the term $v(p)^\beta$
is still well defined as an element of
$(\Q\otimes_\Z\Gamma_v)_\circ$, so the above inequality
makes sense in any case.) Fix $v$, $m$, $\beta$ as in
the corollary, and a distinguished element
$\underline\alpha:=(\alpha_n~|~n\in\N)\in\Ker\,u_A$,
and define the sequence of rational numbers
$(\lambda_{n+1}~|~n\in\N)$ as in
\eqref{subsec_approx-with-Teich}, as well as the subset
$S_n\subset\N[1/p]^{\oplus 2}$ and the $\bar A$-submodule
$\cI(e,n,\gamma)$ of $\Fil^\gamma A_U$, for every
$e\in\Fil^\gamma\bE_U$ and every integer $n>0$. Notice that
\set\begin{equation}\label{eq_lim-of-lambdas}
\lim_{n\to+\infty}\lambda_n=p/(p-1).
\end{equation}
We apply proposition \ref{prop_Scholze-approx-lemma} to
the perfectoid rings with $\Z[1/p]$-graded structure $\cE_U$
and $\cA_U$ given by example \ref{ex_filtered-type} : for
every integer $n>0$ we then get $e_n\in\Fil^\gamma\bE_U$,
$d_n\in\Fil^\gamma A_U$, a finite subset $T_n\subset S_n$, and
for every $(i,j)\in T_n$ an element $b_{i,j,n}\in\Fil^{(1-j)\gamma}A_U$ 
such that
$$
a=\phi^\flat_U(e_n)+\sum_{(i,j)\in T_n}
b_{i,j,n}\cdot\phi^\flat_U(\alpha_0^i\cdot e_n^j)+p^nd_n.
$$
Pick also $\eps\in\N[1/p]\setminus\{0\}$ such that
$\alpha_0\in\Fil^\eps\bE_U$, and $k\in\N$ large enough
so that $\gamma+k\eps>0$.

\begin{claim}\label{cl_of-betas-and-m}
We have $v(p^nd_n)\leq v(p)^{\beta+m}$ for every
sufficiently large $n\in\N$.
\end{claim}
\begin{pfclaim} We have $p^kd_n\in\Fil^{\gamma+k\eps}A_U$,
and notice that $\Fil^\delta A_U\subset A^{\circ\circ}_U$ for
every $\delta>0$. It follows that $v(p^kd_n)<1$, so the
claim holds for every $n\geq k+m+\beta$.
\end{pfclaim}

Recall that $pu=\phi^\flat_U(\alpha_0)$ for some
$u\in\bar A{}^\times$; suppose first that there exists
$n_0\in\N$ such that $v(\phi^\flat_U(e_n))\leq v(p)^m$
for every $n\geq n_0$. Then we have
$$
v(b_{i,j,n}\cdot\phi^\flat_U(\alpha_0^i\cdot e_n^j))\leq
v(b_{i,j,n})\cdot v(u)^i\cdot v(p)^{i+mj}
\qquad
\text{for every $n\geq n_0$ and every $(i,j)\in T_n$}
$$
and on the other hand
$b_{i,j,n}\cdot\phi^\flat_u(\alpha_0^{(1-j)k})\in
\Fil^{(1-j)(\gamma+k\eps)}\subset A_U^{\circ\circ}$, and
$v(u)^i\cdot v(x)<1$ for every $x\in A_U^{\circ\circ}$, whence
$$
v(u)^i\cdot v(b_{i,j,n}\cdot\phi^\flat_u(\alpha_0^{(1-j)k}))<1
\qquad
\text{for every $n\in\N$ and every $(i,j)\in T_n$}.
$$
By the same token, we have $v(u)^{-(1-j)k}\cdot v(p)^\delta<1$ for
every $\delta>0$, whence
$$
v(p)^\delta\leq v(\phi^\flat_U(\alpha_0^{(1-j)k}))
\qquad
\text{for every $\delta>(1-j)k$}.
$$
Summing up we obtain, for every $(i,j)\in T_n$, the inequality
$$
v(b_{i,j,n}\cdot\phi^\flat_U(\alpha_0^i\cdot e_n^j))\leq v(p)^{\beta+m}
\qquad\text{whenever}\qquad
i+mj-(1-j)k>m+\beta.
$$
The latter condition is the same as the inequality
$i>(1-j)(m+k)+\beta$. However, notice that $1-j<i/n$,
so this inequality holds provided $i>i\cdot(m+k)/n+\beta$,
{\em i.e.} when
\set\begin{equation}\label{eq_push-up}
i\cdot(1-(m+k)/n)>\beta.
\end{equation}
Lastly, from \eqref{eq_lim-of-lambdas} we see that
$\beta<\lambda_n$ for every sufficiently large $n$,
and since $i\geq\lambda_n$, we may find $n$ large
enough so that \eqref{eq_push-up} is verified for
every $(i,j)\in T_n$; taking into account claim
\ref{cl_of-betas-and-m}, we get the corollary, in this case.

Thus, we may assume that there exists an infinite subset
$\Sigma\subset\N$ such that $v(p)^m<v(\phi^\flat_U(e_s))$
for every $s\in\Sigma$, and taking into account claim
\ref{cl_of-betas-and-m}, we are reduced to checking that
$$
v(b_{i,j,s}\cdot\phi^\flat_U(\alpha_0^i\cdot e_s^j))\leq
v(p)^\beta\cdot v(\phi_U^\flat(e_s))
\qquad
\text{for some $s\in\Sigma$ and every $(i,j)\in T_s$}.
$$
However, arguing as in the foregoing we easily see that
$$
v(\phi^\flat_U(\alpha_0^i))\leq v(p)^\delta
\qquad
\text{for every $\delta<i$}.
$$
Hence, after replacing $\beta$ by some $\beta'\in\N[1/p]$
with $p/(p-1)>\beta'>\beta$, we reduce to checking that
there exists $s\in\Sigma$ such that
$$
v(b_{i,j,s})\cdot v(p)^{i-\beta}\leq v(p)^{(1-j)m}
\qquad
\text{for every $(i,j)\in T_s$}.
$$
On the other hand, we have
$b_{i,j,n}\cdot\phi_U^\flat(\alpha_0^{(1-j)k})\in A^{\circ\circ}_U$,
and $v(u)^{(1-j)k}\cdot v(x)<1$ for every $x\in A_U^{\circ\circ}$,
whence
$$
v(b_{i,j,n})\cdot v(p)^{(1-j)k}<1
\qquad
\text{for every $n\in\N$ and every $(i,j)\in T_n$}.
$$
Thus, it suffices to show that there exists $s\in\Sigma$
such that
$$
i-\beta>(1-j)k+(1-j)m=(1-j)(m+k)
\qquad
\text{for every $(i,j)\in T_s$}.
$$
This condition is equivalent to \eqref{eq_push-up},
so we may argue as in the foregoing to conclude.
\end{proof}

\sset\subsubsection{}\label{subsec_alternative-meth}
Now, in the situation of corollary \ref{cor_Scholze-approx-lemma},
let
$$
\textstyle{C(f_\bullet):=\Cont(A_U)\cap
R_{A_U}\big(\frac{f_1}{f_0},\cdots,\frac{f_n}{f_0}\big)}
\qquad\text{and}\qquad
S(f_\bullet):=C(f_\bullet)\cap\Spa\,\underline A
$$
for a sequence $f_\bullet:=(f_0,\dots,f_n)$ of elements of
$A_U$ that generates an open ideal. There exist
$a_0,\dots,a_n\in A_U$ and an integer $m>0$ such that
$p^{m-1}=\sum_{i=0}^na_if_i$, and after replacing each $a_i$
by $p^na_i$ for a suitable $n\in\N$, and $m$ by $m+n$, we
may even assume that $a_0,\dots,a_n\in\bar A$ (recall that
$\bar A$ is the image of $A$ in $A_U$, and $\bar\bE$ is the
image of $\bE$ in $\bE_U$) whence
\set\begin{equation}\label{eq_Naoko-0}
v(p)^m\leq\max(v(pa_i)\cdot v(f_i)~|~i=0,\dots,n)\leq
\max(v(f_i)~|~i=0,\dots,n)
\end{equation}
for every $v\in\Cont(A_U)$. It follows that
\set\begin{equation}\label{eq_Naoko-1}
v(f_0)\geq v(p)^m
\qquad
\text{for every $v\in C(f_\bullet)$}.
\end{equation}
By corollary \ref{cor_Scholze-approx-lemma} we may find
$e_0,\dots,e_n\in\bE_U$ such that for every $i=1,\dots,n$
we have :
\set\begin{equation}\label{eq_Naoko-2}
v(f_i-\phi^\flat_U(e_i))\leq v(p)\cdot\max(v(\phi^\flat_U(e_i)),v(p)^m)
\qquad
\text{for every $v\in\Cont(A_U)$}.
\end{equation}
We may now state the announced generalization of theorem
\ref{th_Scholze-tilt} :

\begin{proposition}\label{prop_Scholze-tilt}
With the notation of \eqref{subsec_alternative-meth}, we have :
\begin{enumerate}
\item
The map $\Cont(\phi^\flat_U)$ of lemma
{\em\ref{lem_tilt-Scholze}(ii)} is a homeomorphism.
\item
The system $e_\bullet:=(e_0,\dots,e_n)$ generates an open
ideal of\/ $\bE_U$.
\item
$C(f_\bullet)=\Cont(\phi^\flat_U)^{-1}C(e_\bullet)$, 
where $C(e_\bullet):=R_{\bE_U}(\frac{e_1}{e_0},\cdots,\frac{e_n}{e_0})
\cap\Cont(\bE_U)$.
\item
$\phi^\flat_U(e_0)\in\cO^\loc_{\Spa\,\underline A}(S(f_\bullet))^\times$
and $f_i/f_0-\phi^\flat_U(e_i)/\phi^\flat_U(e_0)\in
\cO^\loc_{\Spa\,\underline A}(S(f_\bullet))^{\circ\circ}$
for $i=1,\dots,n$ (notation of \eqref{subsec_rational-site}).
\end{enumerate}
\end{proposition}
\begin{proof}(ii): Set
$\phi^\flat_U(e_\bullet):=(\phi^\flat_U(e_0),\dots,\phi_U^\flat(e_n))$;
to begin with, we remark :

\begin{claim}\label{cl_extend-cor-taut-two}
(ii) holds if and only if the system $\phi^\flat_U(e_\bullet)$
generates an open ideal of $A_U$.
\end{claim}
\begin{pfclaim} Let $b_1,\dots,b_r$ be a finite system of
elements of $\bar\bE$ that generates an ideal of definition.
Then we may find $k\in\N$ large enough so that
$b^k_ie_j\in\bar\bE$ for every $i=1,\dots,r$ and every
$j=1,\dots,n$. Clearly the system
$b^k_\bullet e_\bullet:=(b^k_ie_j~|~i\leq r,\ j\leq n)$ generates
an open ideal of $\bE_U$ if and only if the same holds for
$e_\bullet$; likewise; $\phi^\flat(e_\bullet)$ generates an
open ideal of $A_U$ if and only if the same holds for
$\phi^\flat_U(b^k_\bullet e_\bullet)$. So we may replace $e_\bullet$
by $b^k_\bullet e_\bullet$, and assume from start that
$e_i\in\bar\bE$ for $i=1,\dots,n$, in which case
$\phi^\flat_U(e_i)=\bar u_{\bar A}(e_i)\in\bar A$ for
$i=1,\dots,n$. Now, let $I_\bE\subset\bar\bE$ be the ideal
generated by $e_\bullet$; then $I_\bE$ is open in $\bar\bE$ if
and only if $Z:=\Spec\,\bar\bE/I_\bE$ lies in the non-analytic
locus $X$ of $\Spec\,\bar\bE$; likewise, $I_\bE\bE_U$ is open
in $\bE_U$ if and only if $Z':=\Spec\,\bE_U/I_\bE\bE_U$ lies
in the non-analytic locus $X'$ of $\Spec\,\bE_U$ (lemma
\ref{lem_deja-vu}(v)). On the other hand, let
$j:\bar\bE\to\bE_U$ be the inclusion map; then
$Z'=(\Spec\,j)^{-1}Z$, and lemma \ref{lem_deja-vu}(iii)
implies that $Z$ lies in $X$ if and only if $Z'$ lies
in $X'$. Summing up, we see that $I_\bE$ is open in $\bar\bE$
if and only if $I_\bE\bE_U$ is open in $\bE_U$. Likewise, let
$I_A\subset\bar A$ be the ideal generated by
$\phi^\flat_U(e_\bullet)$; then $I_A$ is open in $\bar A$ if
and only if $I_AA_U$ is open in $A_U$. Thus, we are further
reduced to showing that $e_\bullet$ generates an open ideal
of $\bar\bE$ if and only if $\phi^\flat_U(e_\bullet)$ generates
an open ideal of $\bar A$. The latter assertion follows from
corollary \ref{cor_taut-two}(ii).
\end{pfclaim}

By claim \ref{cl_extend-cor-taut-two}, it suffices to show
that the system $\phi^\flat_U(e_\bullet)$ generates an open ideal
of $A_U$. To this aim, we apply the criterion of lemma
\ref{lem_criterion-opennes}, which reduces to checking
that for every $v\in\Cont(A_U)_\mathrm{a}$ there exists
$i\in\{0,\dots,n\}$ such that $v(\phi^\flat_U(e_i))\neq 0$.
However, from \eqref{eq_Naoko-2} we deduce that for every
$i=0,\dots,n$ and every such $v$ we have :
\begin{itemize}
\item
either $v(f_i)=v(\phi^\flat_U(e_i))$
\item
or else $v(p)>0$ and $v(f_i),v(\phi^\flat_U(e_i))<v(p)^m$.
\end{itemize}
Now, since the system $f_\bullet$ generates an open ideal
of $A_U$, lemma \ref{lem_criterion-opennes} implies that
$v(f_i)\neq 0$ for some $i\leq n$. If $v(p)=0$ we deduce
that $v(\phi_U^\flat(e_i))\neq 0$, as required.

Thus, suppose $v(p)>0$; from \eqref{eq_Naoko-0} we know
that $v(f_i)\geq v(p)^m>0$ for some $i\leq n$, whence
$v(\phi^\flat_U(e_i))=v(f_i)\neq 0$, again as needed.

(iii): In view of \eqref{eq_Naoko-1} and the foregoing,
it follows that
\set\begin{equation}\label{eq_crapule}
v(f_0)=v(\phi^\flat_U(e_0))\geq v(p)^m
\qquad
\text{for every $v\in C(f_\bullet)$}.
\end{equation}
Now, let $v\in C(f_\bullet)$ and $i\in\{0,\dots,n\}$; if
$v(f_i)\geq v(p)^m$, the foregoing also implies that
$$
v(\phi^\flat_U(e_i))=v(f_i)\leq v(f_0)=v(\phi^\flat_U(e_0)).
$$
If $v(f_i)<v(p)^m$, likewise we get
$v(\phi^\flat_U(e_i))<v(p)^m\leq v(\phi^\flat_U(e_0))$, so
$v\in\Cont(\phi^\flat_U)^{-1}C(e_\bullet)$.

Lastly, let $v\in\Cont(A_U)\setminus C(f_\bullet)$; if $v(p)=0$,
we know that $v(f_i)=v(\phi^\flat_U(e_i))$ for every $i=0,\dots,n$,
whence $v\notin\Cont(\phi^\flat_U)^{-1}C(e_\bullet)$. Otherwise,
\eqref{eq_Naoko-0} says that $v(f_j)\geq v(p)^m$ for some
$j\leq n$, whence $v(\phi^\flat_U(e_j))\geq v(p)^m$. If now
$v(f_0)<v(p)^m$, we get $v(\phi^\flat_U(e_0))<v(p)^m$, so
$v\notin\Cont(\phi^\flat_U)^{-1}C(e_\bullet)$. Hence, suppose
$v(\phi^\flat_U(e_0))=v(f_0)\geq v(p)^m>0$. Then there exists
$k\leq n$ such that $v(f_k)>v(f_0)$, whence
$v(f_k)=v(\phi^\flat_U(e_k))$, and again we conclude that
$v\notin\Cont(\phi^\flat_U)^{-1}C(e_\bullet)$.

(i): First, let us show that the continuous map
$\Cont(\phi^\flat_U)$ is bijective. Indeed, a simple
inspection yields a commutative diagram of continuous maps
$$
\xymatrix{ \Cont(A_U) \ar[rr]^-{\Cont(\phi^\flat_U)} \ar[d] & &
\Cont(\bE_U) \ar[d] \\
\Cont(\bar A) \ar[rr]^-{\Cont(\bar u_{\bar A})} & & \Cont(\bar\bE)
}$$
where the vertical arrows are induced by the open inclusion
maps $\bar A\to A_U$ and $\bar\bE\to\bE_U$ (notation of
\eqref{subsec_general-case}). Therefore, combining
theorem \ref{th_Scholze-tilt}(iii) and proposition
\ref{prop_Cont-open-subring}(ii) we already see that
$\Cont(\phi^\flat_U)$ restricts to a bijection
$\Cont(A_U)_\mathrm{a}\isom\Cont(\bE_U)_\mathrm{a}$.
On the other hand, proposition \ref{prop_general-case}(iv)
implies easily that $\Cont(\phi^\flat_U)$ restricts as well to
a bijection $\Cont(A_U)_\mathrm{na}\isom\Cont(\bE_U)_\mathrm{na}$,
whence the contention. Lastly, (ii) and (iii) imply easily
that the topology of $\Cont(A_U)$ is induced by the topology
of $\Cont(\bE_U)$, via the map $\Cont(\phi^\flat_U)$; the
assertion follows.

(iv): Since $f_0$ is invertible in
$\cO^\loc_{\Spa\,\underline A}(S(f_\bullet))$, combining
\eqref{eq_crapule} and proposition \ref{prop_crit-invertible}
we see that the same holds for $\phi^\flat_U(e_0)$. For the
second assertion, in view of corollary \ref{cor_corcor}(ii)
it suffices to show :

\begin{claim}
$v(f_i/f_0-\phi^\flat_U(e_i)/\phi^\flat_U(e_0))<1$ for every
$i=1,\dots,n$ and every $v\in C(f_\bullet)$.
\end{claim}
\begin{pfclaim}[] From \eqref{eq_Naoko-2} we easily deduce that
$$
v(1/f_0-1/\phi^\flat_U(e_0))\leq v(p/f_0)
\qquad
\text{for every $v\in C(f_\bullet)$}
$$
whence $v(f_i/f_0-f_i/\phi^\flat_U(e_0))\leq v(pf_i/f_0)<1$
for every $v\in C(f_\bullet)$ and every $i=1,\dots,n$. We
are then reduced to checking that
$v(f_i/\phi^\flat_U(e_0)-\phi^\flat_U(e_i)/\phi^\flat_U(e_0))<1$,
or equivalently, that $v(f_i-\phi^\flat_U(e_i))<v(\phi^\flat_U(e_0))$
for every such $v$ and $i$. Now, if $v(f_i)<v(p)^m$, we know
already that $v(\phi^\flat_U(e_i))<v(p)^m$, whence
$v(f_i-\phi^\flat_U(e_i))<v(p)^m$, and the assertion follows
from \eqref{eq_crapule} in this case. Otherwise, we have
$v(f_i)=v(\phi^\flat_U(e_i))$, whence
$v(\phi^\flat_U(e_i))\leq v(f_0)=v(\phi^\flat_U(e_0))$,
and thus $v(f_i-v(\phi^\flat_U(e_i)))\leq
v(p)\cdot\max(\phi^\flat_U(e_0),v(p)^m)=
v(p\cdot\phi^\flat_U(e_0))<v(\phi^\flat_U(e_0))$.
\end{pfclaim}
\end{proof}

\sset\subsubsection{}\label{subsec_graded-perfectoid}
Let $\Gamma,\Delta$ be two monoids and
$(A,\underline B)$ a topological ring with $\Delta$-graded
structure. Then we may form the complete topological ring
with $\Delta\oplus\Gamma$-graded structure
$$
(A,\underline B)[\Gamma]^\wedge
$$
by combining the constructions of example
\ref{ex_monoid-alg-graded}(i) and remark
\ref{rem_graded-top-algs}(iv). Explicitly, a direct
inspection shows that $(A,\underline B)[\Gamma]^\wedge$
is the pair $(A[\Gamma]^\wedge,\underline B^c[\Gamma])$,
where $\underline B^c$ is the $\Delta$-graded ring
such that $\gr_\delta B^c$ is the topological closure
of $\gr_\delta B$ in the completion $A^\wedge$ of $A$,
and $A[\Gamma]^\wedge$ is the completion of $A[\Gamma]$,
for the unique topology on the latter ring such that
the inclusion map $A\to A[\Gamma]$ is adic. Notice that
if $A$ is complete and separated, then $B^c=B$, hence
in this case $(A,\underline B)[\Gamma]^\wedge=
(A[\Gamma]^\wedge,\underline B[\Gamma])$.

\begin{lemma}\label{lem_graded-perfectoid}
In the situation of \eqref{subsec_graded-perfectoid},
suppose furthermore that $A$ is perfectoid, and $\Gamma$,
$\Delta$ are both $p$-perfect, and let
$$
A\xrightarrow{\ j_A\ } A[\Gamma]^\wedge\xleftarrow{\ i_A\ }\Gamma
\qquad\text{and}\qquad
\bE(A)\xrightarrow{\ j_\bE\ }\bE(A)[\Gamma]^\wedge
\xleftarrow{\ i_\bE\ }\Gamma
$$
be the natural inclusion maps. Then we have :
\begin{enumerate}
\item
$A[\Gamma]^\wedge$ is perfectoid.
\item
There exists a natural isomorphism of topological
rings with $\Delta\oplus\Gamma$-graded structures
$$
\omega:\bE(A,\underline B)[\Gamma]^\wedge\isom
\bE((A,\underline B)[\Gamma]^\wedge)
$$
fitting into a commutative diagram
$$
\xymatrix{
\bE(A) \ar[rr]^-{j_\bE} \ar[rrd]_{\bE(j_A)} & &
\bE(A)[\Gamma]^\wedge \ar[d]^\omega & &
\Gamma \ar[lld]^{\bE(i_A)} \ar[ll]_-{i_\bE} \\
& & \bE(A[\Gamma]^\wedge).
}$$
\end{enumerate}
\end{lemma}
\begin{proof} (Here $\Gamma$ is identified with the source
$\bE(\Gamma)$ of the morphism $\bE(i_A)$, via $\bar u_\Gamma$.)

(i): Indeed, by remark \ref{rem_p-can-lie-deep}
we may find an ideal of definition $I\subset A$ fulfilling
the conditions of \eqref{subsec_setup-critperf}, and we
let $J:=I^{(p)}$; by proposition
\ref{prop_begin-criterion-perf}(ii), the corresponding
map $\Phi_I:\gr_I^\bullet A\to\gr^\bullet_JA$ is a ring
isomorphism. Set $I':=IA[\Gamma]^\wedge$ and
$J':=JA[\Gamma]^\wedge$; by construction, the topology
of $A[\Gamma]^\wedge$ agrees with the $I'$-adic topology,
and the natural map $(A/I^2)[\Gamma]\to A[\Gamma]^\wedge/I'^2$
is an isomorphism (remark \ref{rem_completion-of-topring}(ii)
and lemma \ref{lem_still-c-adic}(iv)).
Since $\Gamma$ is $p$-perfect, it follows easily
that the Frobenius endomorphism of $A[\Gamma]^\wedge/I'^2$
is surjective, so $A[\Gamma]^\wedge$ is a P-ring.
On the other hand, notice the natural identifications :
$$
\gr^\bullet_{I'}(A[\Gamma]^\wedge)\isom
(\gr^\bullet_IA)[\Gamma]
\qquad
\gr^\bullet_{J'}(A[\Gamma]^\wedge)\isom
(\gr^\bullet_JA)[\Gamma]
$$
which identify the ring homomorphism $\Phi_{I'}$ of
proposition \ref{prop_begin-criterion-perf}(i) with
the map induced by $\Phi_I$ and the $p$-Frobenius
endomorphism of $\Gamma$. Especially, $\Phi_{I'}$ is
an isomorphism, so the assertion follows from theorem
\ref{th_criterium-perfect}.

(ii): According to proposition \ref{prop_graded-perfectoid}(i),
we may find a $\Delta$-graded ideal $I_B$ of $B$ such that
$I_A:=I_BA$ is an ideal of definition of $A$; then notice
that the topological closure $I^\wedge_B$ of $I_B$ in $A$
lies in $I_A$ and $I^\wedge_B\cap B$ is the topological
closure $I^c_B$ of $I_B$ in $B$, therefore $I^\wedge_B=I_A$
and $I^c_B=B\cap I_A$, which says especially that $I^c_B$
is open in $B$. Set $A_0:=A/I_A$ and $B_0:=B/I^c_B$;
since $I^c_B$ is a $\Delta$-graded ideal as well (remark
\ref{rem_graded-top-algs}(i)), we may form the quotient of
$(A,\underline B)[\Gamma]^\wedge$ by the open ideal
$I^c_B[\Gamma]$ of $B[\Gamma]$ as in example
\ref{ex_was-rem-graded-v}(), and get a natural isomorphism :
\set\begin{equation}\label{eq_one-more-graded}
(A,\underline B)[\Gamma]^\wedge\otimes_BB_0\isom
(A_0[\Gamma],B_0[\Gamma]).
\end{equation}
However, the projections
$(A,\underline B)\to(A_0,\underline B_0)$ and
$(A,\underline B)[\Gamma]^\wedge\to
(A,\underline B)[\Gamma]^\wedge\otimes_BB_0$
induce morphisms  (see \eqref{subsec_E-for-not-F_p-alg})
\set\begin{equation}\label{eq_tired-of-grades}
\bE(A,\underline B)\to\bE(A_0,\underline B_0)
\qquad
\bE((A,\underline B)[\Gamma]^\wedge)\to
\bE((A,\underline B)[\Gamma]^\wedge\otimes_BB_0).
\end{equation}
On the other hand, the projections $A\to A_0$ and
$A[\Gamma]^\wedge\to A_0[\Gamma]$ induce isomorphisms
$$
\bE(A)\isom\bE_0:=\bE(A_0)
\qquad\text{and}\qquad
\bE(A[\Gamma]^\wedge)\isom\bE(A_0[\Gamma])
$$
of topological rings (see \eqref{sec_def-A-tilde}).
Denote
$$
A_0\xrightarrow{\ j_{A_0}\ } A_0[\Gamma]^\wedge
\xleftarrow{\ i_{A_0}\ }\Gamma
\qquad\text{and}\qquad
\bE_0\xrightarrow{\ j_{\bE_0}\ }\bE_0[\Gamma]^\wedge
\xleftarrow{\ i_{\bE_0}\ }\Gamma
$$
the natural inclusion maps; taking into account example
\ref{ex_monoid-alg-graded}(iii), it follows that
\eqref{eq_tired-of-grades} are isomorphisms of topological
rings with graded structures, and in view of
\eqref{eq_one-more-graded} it suffices to exhibit an
isomorphism
$\omega_0:(\bE_0[\Gamma],\bE(\underline B_0)[\Gamma])^\wedge
\isom\bE(A_0[\Gamma],B_0[\Gamma])$ such that
$$
\omega_0\circ j_{\bE_0}=\bE(j_{A_0})
\qquad\text{and}\qquad
\omega_0\circ i_{\bE_0}=\bE(i_{A_0})
$$
(where again, $\Gamma$ is identified with the source
$\bE(\Gamma)$ of $\bE(i_{A_0})$ via $\bar u_\Gamma$).
To this aim, set
$$
J:=\Ker\,(\bar u_{A_0}[\Gamma]:\bE_0[\Gamma]\to A_0[\Gamma])
$$
and notice that, since $A_0$ is a discrete topological
ring, the topology of $\bE_0[\Gamma]$ agrees with the
$J$-adic topology (remark \ref{rem_topology-of-E}(ii));
moreover, $J$ is finitely generated, since $A$ is
perfectoid (details left to the reader).
By corollary \ref{cor_E-and-completion} and remark
\ref{rem_topology-of-E}(v) we deduce a natural
isomorphism of topological rings
$$
\omega_0:\bE_0[\Gamma]^\wedge\isom\bE(A_0[\Gamma])
\qquad\text{such that}\qquad
\bar u_{A_0[\Gamma]}\circ\omega_0=\bar u_{A_0}[\Gamma]^\wedge
$$
where
$\bar u_{A_0}[\Gamma]^\wedge:\bE_0[\Gamma]^\wedge\to A_0[\Gamma]$
is the completion of $\bar u_{A_0}[\Gamma]$, and a simple
inspection shows that $\omega_0$ is a morphism
$(\bE_0[\Gamma],\bE(\underline B_0)[\Gamma])^\wedge\to
\bE(A_0[\Gamma],\underline B_0[\Gamma])$, which then must be 
an isomorphism, again by example
\ref{ex_monoid-alg-graded}(iii). Let us check that
$\omega_0\circ i_{\bE_0}=\bE(i_{A_0})$. Indeed, we have
$$
\bar u_{A_0[\Gamma]}\circ\omega_0\circ i_{\bE_0}=
\bar u_{A_0}[\Gamma]^\wedge\circ i_{\bE_0}=i_{A_0}=
\bar u_{A_0[\Gamma]}\circ\bE(i_{A_0})
$$
whence the sought identity, by adjunction. A similar
argument proves the remaining sought identity.
\end{proof}

We conclude this section with a sample of other graded
perfectoid rings that are analogous to certain constructions
found in \cite{Scho}, and might be interesting for other
purposes, but shall not be needed in the rest of this
treatise.

\begin{example}\label{ex_subsumes}
(i)\ \
Let $(\Delta,0,+)$ be a monoid, $(A,\underline B)$ a
topological ring with $\Delta$-graded structure, and
consider two morphisms of $\Delta$-graded monoids
$$
\Gamma_2\xleftarrow{\ \phi\ }\Gamma_1\xrightarrow{\ \psi\ }B^*
$$
where $B^*$ is the $\Delta$-graded monoid defined as in
\eqref{subsec_reintroduce-B*}. Hence, $B[\Gamma_2]$ is
a $\Delta\oplus\Gamma_2$-graded $\Z$-algebra, and we may
form the $\Delta$-graded $\Z$-algebra $B[\Gamma_2]_{/\Delta}$
corresponding to the morphism
$$
\eta:\Delta\oplus\Gamma_2\to\Delta
\qquad
(\delta,\gamma)\mapsto\delta+|\gamma|
\qquad
\text{for every $\delta\in\Delta$ and $\gamma\in\Gamma_2$}
$$
where $|\cdot|:\Gamma_2\to\Delta$ is the $\Delta$-grading
of $\Gamma_2$. It follows easily that
$$
I:=\Ker\,(B[\Gamma_2]\to\Gamma_2\otimes_{\Gamma_1}B)
$$
is a $\Delta$-graded ideal of $B[\Gamma_2]_{/\Delta}$,
and by combining the constructions of examples
\ref{ex_was-rem-graded-v}(i) and \ref{ex_monoid-alg-graded}(i,ii)
we may define the topological ring with $\Delta$-graded
structure
$$
\Gamma_2\otimes_{\Gamma_1}(A,\underline B):=
(A,\underline B)[\Gamma_2]_{/\Delta}/I.
$$
Explicitly, this is the pair $(A',\underline B')$,
such that $A'$ (resp. $B'$) is the maximal separated
quotient of $\Gamma_2\otimes_{\Gamma_1}A$ (resp. of
$\Gamma_2\otimes_{\Gamma_1}B$), where the latter is endowed
with the unique linear topology such that the natural map
$$
A\to\Gamma_2\otimes_{\Gamma_1}A
\qquad
\text{(resp. $B\to\Gamma_2\otimes_{\Gamma_1}B$)}
$$
is adic. Then, $\gr_\delta B'$ is the topological closure
of the image of $\gr_\delta(B[\Gamma_2]_{/\Delta})$ in $A'$,
for each $\delta\in\Delta$. By construction, we have a
natural morphism of topological rings with $\Delta$-graded
structures (resp. of monoids) 
$$
j_{(A,\underline B)}:(A,\underline B)\to
\Gamma_2\otimes_{\Gamma_1}(A,\underline B)
\qquad
\text{(resp.\ 
$i_{(A,\underline B)}:\Gamma_2\to\Gamma_2\otimes_{\Gamma_1}A$\ )}.
$$

(ii)\ \
Suppose now that $\Delta$ is $p$-perfect, and the topology
of $A$ is complete, separated and coarser than the $p$-adic
topology. Then the topological ring with $\Delta$-graded
structure $\bE(A,\underline B)$ is well defined as in
\eqref{subsec_E-for-not-F_p-alg}. Moreover, from the
morphisms $\phi$ and $\psi$ we derive two morphisms of
monoids
$$
\bE(\Gamma_2)\xleftarrow{\ \phi_\bE\ }\bE(\Gamma_1)
\xrightarrow{\ \psi_\bE\ }\bE(B^*)
\qquad
\text{where $\phi_\bE:=\bE(\phi)$ and $\psi_\bE:=\bE(\psi)$}
$$
so the considerations of (i) can be repeated on
$\bE(A,\underline B)$, and we get a topological
ring with $\bE(\Delta)$-graded structure
$$
\bE(\Gamma_2)\otimes_{\bE(\Gamma_1)}\bE(A,\underline B).
$$
\end{example}

\begin{proposition}\label{prop_subsumes}
In the situation of example {\em\ref{ex_subsumes}}, suppose
moreover that $A$ is perfectoid, and $\Delta$, $\Gamma_1$,
$\Gamma_2$ are $p$-perfect. We have :
\begin{enumerate}
\item
The completion $(\Gamma_2\otimes_{\Gamma_1}A)^\wedge$
of\/ $\Gamma_2\otimes_{\Gamma_1}A$ is perfectoid.
\item
There exists a natural isomorphism of topological rings
with $\Delta$-graded structures :
$$
\omega:\bE((\Gamma_2\otimes_{\Gamma_1}(A,\underline B))^\wedge)\isom
(\Gamma_2\otimes_{\Gamma_1}\bE(A,\underline B))^\wedge
$$
(notation of \eqref{subsec_E-for-not-F_p-alg}) fitting
into the commutative diagrams
$$
\xymatrix{ \bE(A,\underline B) \ar[rr]^-{\bE(j_{(A,\underline B)})}
\ar[rrd]_{j_{\bE(A,\underline B)}} & &
\bE((\Gamma_2\otimes_{\Gamma_1}(A,\underline B))^\wedge)
\ar[d]^\omega & \Gamma_2 \ar[rr]^-{\bE(i_{(A,\underline B)})}
\ar[rrd]_{i_{\bE(A,\underline B)}} & &
\bE((\Gamma_2\otimes_{\Gamma_1}A)^\wedge) \ar[d]^\omega \\
& & (\Gamma_2\otimes_{\Gamma_1}\bE(A,\underline B))^\wedge &
& & (\Gamma_2\otimes_{\Gamma_1}\bE(A))^\wedge.
}$$
\end{enumerate}
\end{proposition}
\begin{proof} (Here we identify $\Gamma_1$ and $\Gamma_2$
with $\bE(\Gamma_1)$ and $\bE(\Gamma_2)$, via the isomorphisms
$\bar u_{\Gamma_1}$ and $\bar u_{\Gamma_2}$.) Let us form as well
the topological ring with $\Delta$-graded structure
$$
(A',\underline D):=
((A,\underline B)[\Gamma_2]^\wedge)_{/\Delta}.
$$
Explicitly, $A'$ is the completion of $A[\Gamma_2]$, where
the latter is endowed with the unique linear topology such
that the natural map $A\to A[\Gamma_2]$ is adic; then
$\gr_\delta D$ is the topological closure in $A'$
of $\gr_\delta(B[\Gamma_2]_{/\Delta})$, for every
$\delta\in\Delta$. Here, as in example \ref{ex_subsumes}(i),
the $\Delta$-graded $\Z$-algebra $B[\Gamma_2]_{/\Delta}$
is obtained from the $\Delta\oplus\Gamma_2$-graded
$\Z$-algebra $B[\Gamma_2]$ via the morphism $\eta$.
With this notation, from proposition
\ref{prop_replaces-Mat-Th.8.1}(i,v) and a simple
inspection we get a natural identification
\set\begin{equation}\label{eq_last-ingredient}
(\Gamma_2\otimes_{\Gamma_1}(A,\underline B))^\wedge
\isom(A',\underline D)/ID.
\end{equation}
Set $(\bE,\underline B_\bE):=\bE(A,\underline B)$; from
lemma \ref{lem_graded-perfectoid} and proposition
\ref{prop_gamma-graded-E-funct}(iii), we know already
that $A'$ is perfectoid, and there exists
an isomorphism
\set\begin{equation}\label{eq_too-many-isos}
(\bE',\underline D_\bE):=\bE(A',\underline D)\isom
((\bE,\underline B_\bE)[\Gamma_2]^\wedge)_{/\Delta}
\end{equation}
of topological rings with $\Delta$-graded structures. Set
$$
\cI:=\Ker\,(B_\bE[\Gamma_2]\to\Gamma_2\otimes_{\Gamma_1}B_\bE).
$$

\begin{claim}\label{cl_running-on-spot}
(i)\ \
$\Phi_{\bE'}(\cI\bE')=\cI\bE'$.
\begin{enumerate}
\addenu
\item
$\{\cI\bE'\}=(IA')^c$.
\end{enumerate}
\end{claim}
\begin{pfclaim} Let
$\bar\psi_\bE:=\bE(i)\circ\psi_\bE:\Gamma_1\to\bE$ where
$i:B^*\to A$ is the natural map. Clearly, the ideals $IA'$
and $\cI\bE'$ are generated respectively by the systems
$$
\cS_A:=(\phi(\gamma)-\psi(\gamma)~|~\gamma\in\Gamma_1)
\qquad\text{and}\qquad
\cS_\bE:=
(\phi_\bE(\gamma)-\bar\psi_\bE(\gamma)~|~\gamma\in\Gamma_1).
$$
Since $\Gamma_1$ is perfect, assertion (i) follows already,
and for (ii) we come down to showing that the topological
closure of the ideal generated by the system
$\bar u_{A'}(\cS_\bE)$ equals the topological closure of the
ideal generated by $\cS_A$. Now, proposition
\ref{prop_combinatorial} says that each generator
$\bar u_{A'}(\phi_\bE(\gamma)-\bar\psi_\bE(\gamma))$
can be written as a $p$-adically convergent series
$\sum_{n\in\N}p^nc_{\gamma,n}$ where
\set\begin{equation}\label{eq_evaluate-at-zero}
c_{\gamma,0}=
\bar u_{A'}(\phi_\bE(\gamma))-\bar u_{A'}(\bar\psi_\bE(\gamma))
\end{equation}
and for every $n>0$, the term $c_{\gamma,n}$ is a finite
$\Z_p$-linear combination of elements of the form
$$
\bar u_{A'}(\phi_\bE(\gamma)^s\cdot\bar\psi_\bE(\gamma)^{s'})-
\bar u_{A'}(\phi_\bE(\gamma)^{s'}\cdot\bar\psi_\bE(\gamma)^s)
$$
where $s,s'\in p^{-n}\N$ and $s+s'=1$.
However, lemma \ref{lem_graded-perfectoid}(ii) implies that
$$
\bar u_{A'}(\phi_\bE(\gamma))=\phi(\gamma)
\qquad\text{and}\qquad
\bar u_{A'}(\bar\psi_\bE(\gamma))=\psi(\gamma)
\qquad
\text{for every $\gamma\in\Gamma_1$}.
$$
Thus $\{\cI\bE'\}$ lies in the topological closure
of the ideal generated by the elements
$$
t_{\gamma,s}:=\phi(\gamma)^s\cdot\psi(\gamma)^{s'}-
\phi(\gamma)^{s'}\cdot\psi(\gamma)^s
\qquad
\text{for every $\gamma\in\Gamma_1$ and $s,s'\in\N[1/p]$
with $s+s'=1$}.
$$
To ease notation, set $\lambda:=\phi(\gamma)$ and
$\mu:=\psi(\gamma)$; we notice that
$$
(\lambda^s-\mu^s)\cdot(\mu^{1-s}+\lambda^{1-s})-(\lambda-\mu)=
t_{\gamma,s}
$$
and since $\lambda^s-\mu^s,\lambda-\mu\in\cS_A$, we
conclude that $t_{\gamma,s}\in IA'$. By the same token,
from \eqref{eq_evaluate-at-zero} we get
$c_0=\lambda-\mu$, so finally $(IA')^c=\{\cI\bE'\}+p(IA')^c$.
From this, an easy induction yields
$$
(IA')^c=\{\cI\bE'\}+p^n(IA')^c
\qquad
\text{for every $n\in\N$}
$$
whence the claim, since $\{\cI\bE'\}$ is a closed ideal.
\end{pfclaim}

From claim \ref{cl_running-on-spot} and corollary
\ref{cor_perfectoid-quot} we deduce that $A'/(IA')^c$
is perfectoid and there exists an isomorphism of topological
rings $\omega':\bE(A'/(IA')^c)\isom\bE'/(\cI\bE')^c$ such that
$$
\omega'\circ\bE(\pi_{A'})=\pi_{\bE'}
$$
where $\pi_{A'}:A'\to A'/(IA')^c$ and
$\pi_{\bE'}:\bE'\to\bE'/(\cI\bE')^c$ are the projections.
This completes the proof of (i), and we also see that
the projection $(\bE',\underline D_\bE)\to
(\bE',\underline D_\bE)/\cI\underline D_\bE$ factors
through a morphism
\set\begin{equation}\label{eq_to-make-a-long}
\bE((A',\underline D)/ID)\to
(\bE',\underline D_\bE)/\cI\underline D_\bE
\end{equation}
which must be an isomorphism, by example
\ref{ex_monoid-alg-graded}(iii). On the other hand, from
\eqref{eq_too-many-isos} we get a natural identification
\set\begin{equation}\label{eq_story-short}
(\bE',\underline D_\bE)/\cI\underline D_\bE\isom
(\Gamma_2\otimes_{\Gamma_1}\bE(A,\underline B))^\wedge.
\end{equation}
Lastly, combining \eqref{eq_to-make-a-long},
\eqref{eq_story-short} and \eqref{eq_last-ingredient}
we get the sought isomorphism $\omega$. The commutativity
of the diagrams in (ii) follows by direct inspection,
taking into account the corresponding properties of
the isomorphism \eqref{eq_too-many-isos} stated in lemma
\ref{lem_graded-perfectoid}(ii) : details left to the
reader.
\end{proof}

\sset\subsubsection{}\label{subsec_perf-blowup}
Let $\Lambda$ be any (small) set and $(\Delta,+,0,\leq)$
any ordered abelian group; we set
$$
\Delta^+=\{\delta\in\Delta~|~\delta\geq 0\}
\qquad
P:=\Delta^{+(\Lambda)}
\qquad
Q:=\Delta^{+(\Lambda\times\Lambda)}
$$
(cp. definition \ref{def_ordered-group}(i)). Denote by
$I\subset P$ the ideal generated by the canonical system
of generators $e_\bullet:=(e_\lambda~|~\lambda\in\Lambda)$
for $P$, {\em i.e.}
$e_\lambda:=(\delta_{\lambda\mu}~|~\mu\in\Lambda)$,
where $\delta_{\lambda\mu}:=1$ for $\lambda=\mu$, and
$\delta_{\lambda\mu}:=0$ otherwise, for every
$\lambda,\mu\in\Lambda$. Likewise, we let
$$
(f_\lambda,f'_\lambda~|~\lambda\in\Lambda)
\qquad
\text{(resp. $(e_{\lambda,\lambda'}~|~\lambda,\lambda'\in\Lambda)$)}
$$
be the canonical system of generators of $P\times P$
(resp. of $Q$).

$\bullet$\ \
With this notation, we define $\Delta$-gradings on
$P$, $P\times P$ and $Q$ by the rules :
$$
|\delta e_\lambda|_P:=\delta
\qquad
|\delta f_\lambda+\delta'f'_\lambda|_{P\times P}:=\delta
\qquad
|\delta e_{\lambda,\lambda'}|_Q:=\delta
$$
for every $\lambda,\lambda'\in\Lambda$ and every
$\delta,\delta'\in\Delta$ (see definition
\ref{def_grad-monoids}(i)), and we set
$$
I^{\La\delta\Ra}:=\{x\in P~|~|x|_P\geq\delta\}
\qquad
\text{for every $\delta\in\Delta$}.
$$
Notice that, in case $\Delta=\Z[1/p]$, the ideals
$I^{\La\delta\Ra}$ of $P$ are the same as the ones
introduced in \eqref{subsec_mon-fract-powers}, so
the notation does not conflict with {\em loc.cit.}
Lastly, we define
$$
\Gamma^+_{\!I}:=
\bigcup_{\delta\in\Delta^+}I^{\La\delta\Ra}\times\{\delta\}
\qquad
\Gamma_{\!I}:=
\bigcup_{\delta\in\Delta}I^{\La\delta\Ra}\times\{\delta\}
$$
and notice that $\Gamma^+_{\!I}$ (resp. $\Gamma_{\!I}$) is
a submonoid of $P\times\Delta^+$ (resp. of $P\times\Delta$),
so it inherits a natural $\Delta$-grading, given by the
projection on the factor $\Delta$. We have morphisms of
$\Delta$-graded monoids
\set\begin{equation}\label{eq_present-Gamma}
{\diagram
Q \ar@<.5ex>[r]^-\psi \ar@<-.5ex>[r]_-{\psi'} &
P\times P \ar[r]^-\phi & \Gamma^+_{\!I}
\enddiagram}
\end{equation}
such that
$$
\phi(\delta f_\lambda+\delta' f'_\lambda):=
((\delta+\delta')e_\lambda,\delta)
\qquad
\psi(\delta e_{\lambda,\lambda'}):=\delta(f_\lambda+f'_{\lambda'})
\qquad
\psi'(\delta e_{\lambda,\lambda'}):=\delta(f'_\lambda+f_{\lambda'})
$$
for every $\lambda,\lambda'\in\Lambda$ and every
$\delta,\delta'\in\Delta$.

\begin{lemma}\label{lem_present-Gamma}
With the notation of \eqref{subsec_perf-blowup}, we have :
\begin{enumerate}
\item
\eqref{eq_present-Gamma} is a presentation for\/
$\Gamma^+_{\!I}$, {\em i.e.} $\phi$ identifies
$\Gamma^+_{\!I}$ with the coequalizer of $\psi$ and
$\psi'$.
\item
Define a $\Delta$-grading of $P\times\Delta^+$ by the rule :
$(x,\delta)\mapsto|x|_P-\delta$ for every $x\in P$ and
$\delta\in\Delta^+$. Then the mapping
$$
P\times\Delta^+\to\Gamma_{\!I}
\qquad
(x,\delta)\mapsto(x,|x|_P-\delta)
$$
is an isomorphism of $\Delta$-graded monoids.
\end{enumerate}
\end{lemma}
\begin{proof}(i): Indeed, say that
$$
g:=\sum_{\lambda\in\Lambda}(r_\lambda f_\lambda+r'_\lambda f'_\lambda)
\qquad\text{and}\qquad
h:=\sum_{\lambda\in\Lambda}(s_\lambda f_\lambda+s'_\lambda f'_\lambda)
$$
are two elements of $P\times P$ whose images agree
in $\Gamma^+_{\!I}$, so that
$r_\lambda,r'_\lambda,s_\lambda,s'_\lambda\in\Delta^+$ and
\set\begin{equation}\label{eq_present-gamma}
r_\lambda+r'_\lambda=s_\lambda+s'_\lambda
\qquad\text{for every $\lambda\in\Lambda$}
\qquad\text{and}\qquad
\sum_{\lambda\in\Lambda}r_\lambda=
\sum_{\lambda\in\Lambda}s_\lambda.
\end{equation}
Let $\sim$ be the equivalence relation defined by the
pair $(\psi,\psi')$; we have to check that $g\sim h$.
Now, consider any $\lambda\in\Lambda$, and suppose
that $r_\lambda\leq s_\lambda$; then we may write
$g=g'+r_\lambda f_\lambda$, $h=h'+r_\lambda f_\lambda$
for elements $g',h'\in P\times P$ such that
$\phi(g')=\phi(h')$; it then suffices to check that
$g'\sim h'$, and we are reduced to the case where
$r_\lambda=0$. In case $r_\lambda>s_\lambda$, we argue
symmetrically, to reduce to the case where $s_\lambda=0$.
Likewise, we can reduce to the case where either
$r'_\lambda=0$ or $s'_\lambda=0$.
Repeating the argument for every $\lambda\in\Lambda$,
and taking \eqref{eq_present-gamma} into account,
we may suppose from start that
$$
g=\sum_{\lambda\in A}r_\lambda f_\lambda+\sum_{\mu\in B}r_\mu f'_\mu
\qquad
h=\sum_{\lambda\in A}r_\lambda f'_\lambda+\sum_{\mu\in B}r_\mu f_\mu
$$
for some finite subsets $A,B\in\Lambda$ with
$A\cap B=\emptyset$, and a system
$(r_\lambda~|~\lambda\in A\cup B)$ of elements of
$\Delta^+\setminus\{0\}$ such that
\set\begin{equation}\label{eq_combinat-on-high}
\sum_{\mu\in B}r_\mu=\sum_{\lambda\in A}r_\lambda.
\end{equation}
Next, we argue by induction on the cardinality $c$ of
$A\cup B$. If $c=0$, there is nothing to prove. Suppose
then that $c>0$ and that the assertion is already known
whenever the cardinality of $A\cup B$ is strictly
smaller than $c$. Notice that $A\neq\emptyset$, since
otherwise \eqref{eq_combinat-on-high} would imply that
$B=\emptyset$, contradicting the assumption that $c>0$;
likewise, we must have $B\neq\emptyset$.
Thus, pick any $\lambda\in A$ and $\mu\in B$, set
$A':=A\setminus\{\lambda\}$, $B':=B\setminus\{\mu\}$,
and suppose first that $r_\lambda\leq r_\mu$; we set
$$
g':=\sum_{\lambda\in A'}r_\lambda f_\lambda+
\sum_{\mu\in B'}r_\mu f'_\mu+(r_\mu-r_\lambda)\cdot f'_\mu
\qquad
h':=\sum_{\lambda\in A'}r_\lambda f'_\lambda+
\sum_{\mu\in B'}r_\mu f_\mu+(r_\mu-r_\lambda)\cdot f_\mu
$$
so that $g=g'+r_\lambda\cdot(f_\lambda+f'_\mu)$
and $h=h'+r_\lambda\cdot(f'_\lambda+f_\mu)$. However,
$r_\lambda\cdot(f_\lambda+f'_\mu)\sim
r_\lambda\cdot(f'_\lambda+f_\mu)$, so we are reduced
to checking that $g'\sim h'$. The latter is known
by our inductive assumption. A symmetric argument
will do, in case $r_\lambda>r_\mu$.

(ii) shall be left to the reader.
\end{proof}

\sset\subsubsection{}\label{subsec_blow-again}
In the situation of \eqref{subsec_perf-blowup}, take
$\Delta:=\Z[1/p]$ with its standard ordering, so that
$\Delta^+=\N[1/p]$. In this case, it is easily seen that
$P$, $\Gamma^+_{\!I}$ and $\Gamma_{\!I}$ are $p$-perfect monoids.
In this paragraph, it will be convenient to switch to a
multiplicative notation for the composition laws of these
monoids : in other words, we shall hereafter replace $(P,+,0)$
by $(\exp P,\cdot,1)$, and likewise for $\Gamma^+_{\!I}$ and
$\Gamma_{\!I}$ (notation of \eqref{sec_toric}). Correspondingly,
the product $\delta\cdot x$ of an element $\delta\in\Delta^+$
and an element $x\in P$ shall be replaced by the exponential
$x^\delta$.

Let also $A$ be any perfectoid ring, set $\bE:=\bE(A)$
and denote by $\underline A$ the topological ring with
$\Delta$-graded structure $(A,A,\Delta)$, such that
$\gr_0A=A$ and $\gr_\delta A=0$ for every $\delta\neq 0$
in $\Delta$; define likewise the topological ring with
$\Delta$-graded structure $\underline\bE:=(\bE,\bE,\Delta)$.
Choose an ideal of definition $J$ for $A$, and set
$\cJ:=\bar u_{A/pA}^{\ -1}(J/pA)$, which is an ideal
of definition for $\bE$.
Let $\beta_\bullet:=(\beta_\lambda~|~\lambda\in\Lambda)$
be any system of elements of $\bE$; then $\beta_\bullet$
determines a unique morphism of monoids
$$
\phi_\bE:P\to\bE
\qquad\text{such that}\qquad
e^\delta_\lambda\mapsto\beta^\delta_\lambda
\qquad
\text{for every $\lambda\in\Lambda$ and $\delta\in\Delta$}
$$
and composing with the map $\bar u_A:\bE\to A$ we
get also a morphism $\phi_A:P\to A$. Moreover, notice
that $I^{\La 0\Ra}=P$, so we have as well a morphism
$i_0:P\to\Gamma^+_{\!I}$, and we can define the rings
$$
\cE_+:=\Gamma^+_{\!I}\otimes_P\bE\subset
\cE:=\Gamma_{\!I}\otimes_P\bE
\qquad\text{and}\qquad
\cA_+:=\Gamma^+_{\!I}\otimes_PA\subset
\cA:=\Gamma_{\!I}\otimes_PA
$$
as well as the topological rings with $\Delta$-graded
structures
$$
\begin{aligned}
(\cE^\wedge_+,\underline\cE^c_+):=\, &
(\Gamma^+_{\!I}\otimes_P\underline\bE)^\wedge
& \subset & &
(\cE^\wedge,\underline\cE^c_+):=
(\Gamma_{\!I}\otimes_P\underline\bE)^\wedge \\
(\cA^\wedge_+,\underline\cA^c_+):=\, &
(\Gamma^+_{\!I}\otimes_P\underline A)^\wedge
& \subset & &
(\cA,\underline\cA^c):=
(\Gamma_{\!I}\otimes_P\underline A)^\wedge
\end{aligned}
$$
(notation of \eqref{subsec_restrict-scalar-mon} and
example \ref{ex_subsumes}). We endow $\cE$ and $\cE_+$
with their $\cJ$-adic topologies, and $\cA$, $\cA_+$
with their $J$-adic topologies; then $\cE^\wedge_+$ (resp.
$\cE^\wedge$) is the completion of $\cE_+$ (resp. of $\cE$)
and $\cA^\wedge_+$ (resp. $\cA^\wedge$) is the completion
of $\cA_+$ (resp. of $\cA$).
Furthermore, the rings $\cA$ and $\cE$ inherit natural
$\Delta$-gradings from $\Gamma_{\!I}$, such that
$$
I^{\La\delta\Ra}\otimes_P\bE=\gr_\delta\cE
\qquad
I^{\La\delta\Ra}\otimes_PA=\gr_\delta\cA
\qquad
\text{for every $\delta\in\Delta$}
$$
and by restricting to $\Delta^+$, we get corresponding
$\Delta^+$-gradings on $\cA_+$ and $\cE_+$. Then, a
direct inspection shows that $\gr_\delta\cE^c$ (resp.
$\gr_\delta\cA^c$) is the separated completion of
$\gr_\delta\cE$ (resp. of $\gr_\delta\cA$), for every
$\delta\in\Delta$. Let $i_\Gamma:\Gamma^+_I\to\Gamma_I$
be the inclusion map; as in example \ref{ex_subsumes}(i),
we point out the commutative diagram of monoids and of
topological rings with $\Delta$-graded structures :
$$
\xymatrix{
\underline\bE \ar[rr]^-{j^+_\bE} \ar[rrd]_{j_\bE} & &
(\cE^\wedge_+,\underline\cE^c_+) \ar[d]|{(i_\Gamma\otimes_P\bE)^\wedge} &
\Gamma^+_I \ar[r]^-{i^+_A} \ar[d]^{i_\Gamma} \ar[l]_-{i^+_\bE} &
(\cA^\wedge_+,\underline\cA^c_+) \ar[d]|{(i_\Gamma\otimes_PA)^\wedge}
& & \underline A \ar[ll]_-{j^+_A} \ar[lld]^{j_A} \\
& & (\cE^\wedge,\underline\cE^c) &
\Gamma_I \ar[r]^-{i_A} \ar[l]_-{i_\bE} &
(\cA^\wedge,\underline\cA^c).
}$$

\begin{remark} In the situation of \eqref{subsec_blow-again},
suppose furthermore that the system $\beta_\bullet$ generates
an open ideal of $\bE$. Then, combining corollaries
\ref{cor_pescato}(ii) and \ref{cor_taut-two}(ii) together
with proposition \ref{prop_a-unif-for-perf}(ii), we deduce that
$\Gamma_{\!I}^+\otimes_P\bE$ and $\Gamma_{\!I}\otimes_P\bE$ are
already $\cJ$-adically complete and separated, and
$\Gamma_{\!I}^+\otimes_PA$ and $\Gamma_{\!I}\otimes_PA$ are
$J$-adically complete and separated.
\end{remark}

\begin{corollary}
With the notation of \eqref{subsec_blow-again}, the
rings $\cE^\wedge_+$, $\cE^\wedge$, $\cA^\wedge_+$ and
$\cA^\wedge$ are perfectoid, and we have natural
isomorphisms of topological rings with $\Delta$-graded
structures
$$
\omega_+:\bE(\cA^\wedge_+,\underline\cA^c_+)
\isom(\cE^\wedge_+,\underline\cE^c_+)
\qquad
\omega:\bE(\cA^\wedge,\underline\cA^c)
\isom(\cE^\wedge,\underline\cE^c)
$$
fitting into the commutative diagrams
$$
\xymatrix{ \underline\bE \ar[rr]^-{\bE(j^+_A)} \ar[rrd]_{j^+_\bE}
& & \bE(\cA^\wedge_+,\underline\cA^c_+) \ar[d]_{\omega_+} & &
\Gamma_I^+ \ar[ll]_-{\bE(i^+_A)} \ar[lld]^-{i^+_\bE} &
\underline\bE \ar[rr]^-{\bE(j_A)} \ar[rrd]_{j_\bE} & &
\bE(\cA^\wedge,\underline\cA^c) \ar[d]_\omega & &
\Gamma_I \ar[ll]_-{\bE(i_A)} \ar[lld]^{i_\bE} \\
& & (\cE^\wedge_+,\underline\cE^c_+) & & & & &
(\cE^\wedge,\underline\cE^c).
}$$
\end{corollary}
\begin{proof} Clearly we have $\bE(\underline A)=\underline\bE$.
The corollary then becomes a special case of proposition
\ref{prop_subsumes}, after we have shown

\begin{claim}
$\phi_\bE\circ\bar u_P=\bE(\bar u_A\circ\phi_\bE)$.
\end{claim}
\begin{pfclaim}[] Clearly
$\phi_\bE\circ\bar u_P=\bar u_\bE\circ\bE(\phi_\bE)$, so we
are reduced to checking the identity $\bar u_\bE=\bE(\bar u_A)$,
and the latter holds by remark \ref{rem_general}(ii).
\end{pfclaim}
\end{proof}

\sset\subsubsection{}
From lemma \ref{lem_present-Gamma}(i) we deduce
a corresponding coequalizer presentation for $\cE_+$
in the category of $\Delta$-graded $\Z$-algebras :
$$
\xymatrix{ \Z[Q] \ar@<.5ex>[r] \ar@<-.5ex>[r] &
\Z[P\times P]\otimes_{\Z[P]}\bE=\bE[P] \ar[r] &
\Z[\Gamma^+_{\!I}]\otimes_{\Z[P]}\bE=\cE_+
}$$
where $\bE$ and $\Z[P\times P]$ are regarded as $\Z[P]$-algebras
via $\phi_\bE$ and respectively via the morphism of monoids
$P\to P\times P$ given by the rule
$e^\delta_\lambda\mapsto f'^\delta_\lambda$ for every
$\lambda\in\Lambda$ and every $\delta\in\Delta^+$.
Likewise we may present $\cA_+$, and there follow
natural isomorphisms
$$
\cE_+\isom\bE[P]/\cI_\bE
\qquad
\cA_+\isom A[P]/\cI_A
$$
where $\cI_\bE\subset\bE[P]$ is the ideal generated by
the system
$$
(\beta^r_\lambda\cdot e_\mu^r-\beta^r_\mu\cdot e^r_\lambda
~|~\lambda,\mu\in\Lambda;\ r\in\N[1/p]).
$$
and $\cI_A\subset A[P]$ is the ideal generated by
the system
$$
(\bar u_A(\beta^r_\lambda)\cdot e_\mu^r-
\bar u_A(\beta^r_\mu)\cdot e^r_\lambda
~|~\lambda,\mu\in\Lambda;\ r\in\N[1/p]).
$$
Especially, notice that $\cI_\bE$ (resp. $\cI_A$) is a
graded ideal of the $\Delta$-graded ring $\bE[P]_{/\Delta}$
(resp. $A[P]_{/\Delta}$), arising via the grading
$|\cdot|_P:P\to\Delta$, and
$$
\gr_\delta\cE_+=
(\gr_\delta\bE[P]_{/\Delta})/\gr_\delta\cI_\bE=
\bE[\gr_\delta P]/\gr_\delta\cI_\bE
\qquad
\text{for every $\delta\in\Delta$}.
$$
We can then form the topological rings with
$\Delta$-graded structure
$$
(\cE^s_+,\underline\cE^s_+):=
(\underline\bE[P]_{/\Delta})/\cI_\bE
\qquad
(\cA^s_+,\underline\cA^s_+):=
(\underline A[P]_{/\Delta})/\cI_A
$$
that are the maximal separated quotients $\cE^s_+$ of $\cE_+$
and $\cA^s_+$ of $\cA_+$, endowed with the $\Delta$-graded
$\Z$-algebra structures
$$
\underline\cE^s_+:=
\bigoplus_{\delta\in\Delta^+}(I^{\La\delta\Ra}\otimes_P\bE)^s
\qquad\text{and}\qquad
\underline\cA^s_+:=
\bigoplus_{\delta\in\Delta^+}(I^{\La\delta\Ra}\otimes_PA)^s
$$
where $(I^{\La\delta\Ra}\otimes_P\bE)^s$ denotes the maximal
separated quotient of $I^{\La\delta\Ra}\otimes_P\bE$, for
every $\delta\in\Delta$, and likewise for
$(I^{\La\delta\Ra}\otimes_PA)^s$. It is now straightforward
that the completion maps $\cE^s_+\to\cE_+^\wedge$ and
$\cA^s_+\to\cA_+^\wedge$ factor through natural isomorphisms
$$
(\cE^s_+,\underline\cE^s_+)^\wedge\isom
(\cE^\wedge_+,\underline\cE^c_+)
\qquad\text{and}\qquad
(\cA^s_+,\underline\cA^s_+)^\wedge\isom
(\cA^\wedge_+,\underline\cA^c_+).
$$

\subsection{Perfectoid spaces}
Some basic verifications concerning perfectoid spaces will
employ a certain amount of formal algebraic geometry,
that we collect hereafter. 

\sset\subsubsection{}\label{subsec_proj-with-Q-greaded}
Quite generally, let $R$ be any $\Q_+$-graded ring. We
associate with $R$ a projective scheme
$$
\Proj\,R
$$
as follows. For every $\gamma\in\Q_{>0}$, let
$R^{(\gamma)}:=\bigoplus_{k\in\N}R_{k\gamma}$, which we regard
as an $\N$-graded ring, whose $k$-graded direct summand
is $R_{k\gamma}$, for every $k\in\N$.  Set
$Y^{(\gamma)}:=\Proj\,R^{(\gamma)}$ for every such $\gamma$.
If $n>0$ is any integer, the discussion of
\eqref{subsec_Omegas} yields a natural isomorphism
$$
Y^{(\gamma)}\xrightarrow{\ \Omega_n^{(\gamma)}\ }Y^{(n\gamma)}
$$
of $R_0$-schemes, and it is easily seen that
$$
\Omega^{(n\gamma)}_p\circ\Omega^{(\gamma)}_n=\Omega^{(\gamma)}_{pn}
\qquad
\text{for every integers $n,p>0$}.
$$
Thus, the rule $\gamma\mapsto Y^{(\gamma)}$ yields a well
defined filtered system of isomorphisms, and we let
$\Proj\,R$ be any choice of an $R_0$-scheme representing
the colimit of this system; {\em e.g.}
$\Proj\,R:=Y^{(\gamma_0)}$, for some fixed $\gamma_0\in\Q_{>0}$,
together with the cocone $(\Omega^{(\gamma)}:Y^{(\gamma)}\to
\Proj\,R~|~\gamma\in\Q_{>0})$ resulting from the filtered
system of morphisms $\Omega^{(\bullet)}_\bullet$.

\sset\subsubsection{}
To ease notation, we shall let $Y:=\Proj\,R$. For every
$\gamma\in\Q_{>0}$, every $n\in\N$, and every $f\in R_{n\gamma}$
we have the open subset $D_+(f)\subset Y^{(\gamma)}$ defined as
in \eqref{subsec_projective-spectra}, and we denote also
$D_+(f)\subset Y$ the image of this open subset under
$\Omega^{(\gamma)}$. It is easily seen that the resulting
open subset of $Y$ is independent of the choice of
$\gamma$. We also set $\cO_Y(0):=\cO_Y$ and
$$
\cO_Y(\gamma):=\Omega^{(|\gamma|)}_*\cO_{Y^{(|\gamma|)}}(\gamma/|\gamma|)
\qquad
\text{for every $\gamma\in\Q\setminus\{0\}$}
$$
and we notice that -- according to \eqref{subsec_qcoh-mod-grad} --
the quasi-coherent $\cO_Y$-module $\cO_Y(\gamma)$ restricts to an
invertible $\cO_Y$-module on the open subset
$$
U_\gamma(R):=\bigcup_{f\in R_{|\gamma|}}D_+(f)
\qquad
\text{for every $\gamma\in\Q\setminus\{0\}$}.
$$

\begin{remark} Let $\Gamma\subset\Q_+$ be any submonoid,
and suppose that $R=\Q_+\times_\Gamma R'$ for a
$\Gamma$-graded ring $R'$. Then it is easily seen that
$\cO_Y(\gamma)=0$ and $U_\gamma(R)=\emptyset$ for every
$\gamma\in\Q\setminus\Gamma^\gp$.
\end{remark}

\sset\subsubsection{}
In view of the isomorphisms \eqref{eq_twists}, we get
natural identifications
\set\begin{equation}\label{eq_divide-and-twist}
\cO_Y(n\gamma)\isom\Omega^{(\gamma)}_*\cO_{Y^{(\gamma)}}(n)
\qquad
\text{for every $n\in\Z$ and every $\gamma\in\Q_{>0}$}.
\end{equation}
For every $\gamma,\gamma'\in\Q$ there exists as well a natural
morphism
\set\begin{equation}\label{eq_if-divisor}
\cO_Y(\gamma)\otimes_{\cO_Y}\cO_Y(\gamma')\to\cO_Y(\gamma+\gamma')
\end{equation}
that generalizes \eqref{eq_annother-onne}. Namely, pick a
common divisor, {\em i.e.} a rational number $\delta\in\Q_{>0}$
and integers $p,q\in\N$ such that $p\delta=\gamma$ and
$q\delta=\gamma'$; we have first a natural identification
$$
\cO_Y(\gamma)\otimes_{\cO_Y}\cO_Y(\gamma')\isom
\Omega^{(\delta)}_*\cO_{Y^{(\delta)}}(p)\otimes_{\cO_Y}
\Omega^{(\delta)}_*\cO_{Y^{(\delta)}}(q)\isom\Omega^{(\delta)}_*
(\cO_{Y^{(\delta)}}(p)\otimes_{\cO_Y^{(\delta)}}\cO_{Y^{(\delta)}}(q))
$$
which we combine with the map
$$
\Omega^{(\delta)}_*
(\cO_{Y^{(\delta)}}(p)\otimes_{\cO_Y^{(\delta)}}\cO_{Y^{(\delta)}}(q))\to
\Omega^{(\delta)}_*\cO_{Y^{(\delta)}}(p+q)\isom\cO_Y(\gamma+\gamma')
$$
deduced from \eqref{eq_annother-onne} and the inverse of the
isomorphism \eqref{eq_divide-and-twist}, to obtain
\eqref{eq_if-divisor}, and it is easily seen that the resulting
map is independent of the choice of $\delta$ : details left to
the reader. Especially, \eqref{eq_if-divisor} restricts to an
isomorphism on $U_{\mathrm{gcd}(\gamma,\gamma')}(R)$, where
$\mathrm{gcd}(\gamma,\gamma')$ is the greatest common
divisor of $\gamma$ and $\gamma'$ (here we let
$\mathrm{gcd}(0,0):=0$; actually, it is not difficult to
show that \eqref{eq_if-divisor} is an isomorphism on the
open subset $U_\gamma(R)\cup U_{\gamma'}(R)$). Furthermore, let
$$
\pi^{(\gamma)}:Y^{(\gamma)}\to\Spec\,R_0
\qquad\text{and}\qquad
\pi:Y\to\Spec\,R_0
$$
be the structure morphisms; the morphism \eqref{eq_force-again}
induces a natural map
$$
\pi^*R_\gamma^\sim\isom\Omega^{(\gamma)}_*\pi^{(\gamma)*}R_1^{(\gamma)}
\to\Omega^{(\gamma)}_*\cO_{Y^{(\gamma)}}(1)=\cO_Y(\gamma)
\qquad
\text{for every $\gamma\in\Gamma\setminus\{0\}$}
$$
which restricts to an epimorphism on $U_\gamma(R)$.

\sset\subsubsection{}\label{subsec_functor-Q-graded}
Let $R'$ be another $\Q_+$-graded ring, and $\phi:R\to R'$
a morphism of $\Q_+$-graded rings; set $Y':=\Proj\,R'$,
$Y'^{(\gamma)}:=\Proj\,R'^{(\gamma)}$, and denote by
$\Omega'^{(\gamma)}:Y^{(\gamma)}\to Y'$ and
$\Omega'^{(\gamma)}_n:Y^{(\gamma)}\to Y'^{(n\gamma)}$ the induced
morphism, for every $\gamma\in\Q_{>0}$ and every integer
$n>0$. The restriction of $\phi$
$$
\phi^{(\gamma)}:R^{(\gamma)}\to R'^{(\gamma)}
\qquad
\text{for every $\gamma\in\Q_{>0}$}
$$
is a morphism of $\N$-graded rings, whence a morphism
$\Proj\,\phi^{(\gamma)}:G(\phi^{(\gamma)})\to Y^{(\gamma)}$
(notation of \eqref{subsec_in-the-situat}). Moreover,
it is easily seen that
$$
\Omega'^{(\gamma)}_nG(\phi^{(\gamma)})=G(\phi^{(n\gamma)})
\qquad
\text{for every $\gamma\in\Q_{>0}$ and every integer $n>0$}.
$$
Furthermore, the resulting diagram of schemes
$$
\xymatrix{ G(\phi^{(\gamma)})
\ar[rr]^-{\Proj\,\phi^{(\gamma)}} \ar[d]_{\Omega'^{(\gamma)}_n} & &
Y^{(\gamma)} \ar[d]^{\Omega^{(\gamma)}_n} \\
G(\phi^{(n\gamma)}) \ar[rr]^-{\Proj\,\phi^{(n\gamma)}} & &
Y^{(n\gamma)}
}$$
commutes for every such $\gamma$ and $n$. Thus, we may set
$$
G(\phi):=\Omega^{(\gamma)}(G(\phi^{(\gamma)}))
\qquad
\text{for any $\gamma\in\Q_{>0}$}
$$
and the colimit of the system
$(\Proj\,\phi^{(\gamma)}~|~\gamma\in\Q_+)$ is a well defined
morphism
$$
\Proj\,\phi:G(\phi)\to Y.
$$

\begin{remark}\label{rem_likewise}
(i)\ \
Set $R_+:=\bigoplus_{\gamma\in\Q_{>0}}R_\gamma$, and define
likewise $R'_+$. Let also $f'\in R'_+$ be any homogeneous
element; we claim that $D_+(f')\subset G(\phi)$ if and
only if $f'$ lies in the radical $J$ of the ideal of $R'$
generated by $\phi(R_+)$. Indeed, say that $f'\in R'_\gamma$
for some $\gamma\in\Q_{>0}$; if $D_+(f')\subset G(\phi)$,
then $D_+(f')\subset G(\phi^{(\gamma)})$ as well, so $f'$
lies in the radical of the ideal of $R'^{(\gamma)}$ generated
by $\phi^{(\gamma)}(R^{(\gamma)}_+)$ (see the discussion of
\eqref{subsec_in-the-situat}), and therefore $f'$ lies
in $J$ as well. Conversely, if $f'\in J$, we may find
$n,k\in\N$, and elements $\nu(i)\in\Q_{>0}$, $a_i\in R_{\nu(i)}$,
$f'_i\in R'_{\gamma-\nu(i)}$ for $i=1,\dots,k$, such that
$f'^n=\sum_{i=1}^kf_i\cdot\phi(a_i)$. Then, pick any common
divisor $\delta$ of $\gamma,\nu(1),\dots,\nu(k)$ in
$\Q_{>0}$; it follows that $f'$ lies in the radical of
the ideal of $R'^{(\delta)}$ generated by
$\phi^{(\delta)}(R_+^{(\delta)})$, whence
$D_+(f')\subset G(\phi^{(\delta)})$, and finally
$D_+(f')\in G(\phi)$. 

(ii)\ \
Especially, if $\phi(R_+)$ generates the ideal $R'_+$
of $R'$, then $G(\phi)=Y'$.
\end{remark}

\sset\subsubsection{}\label{subsec_angular-blow-up}
We shall apply the foregoing constructions to the situation
contemplated in \eqref{subsec_link-with-a-un-etc} : we let
$(A,u_0)$ be a pair consisting of a ring $A$, and a ring
homomorphism
$$
u_0:R_{r,0}:=\Z[T^{1/p^\infty}_1,\dots,T^{1/p^\infty}_r]\to A.
$$
Also, we denote by $T_\bullet P_r\subset P_r:=\N[1/p]^{\oplus r}$
the ideal generated by $\N^{\oplus r}$, we set $\Gamma:=\N[1/p]$,
and we consider the $\Gamma$-graded {\em angular Rees algebra}
$$
\sR\La A\Ra:=\bigoplus_{\gamma\in\Gamma}T_\bullet^{\La\gamma\Ra}A.
$$
Hence, for every $\gamma\in\Gamma$, the direct summand
$\sR\La A\Ra_\gamma$ is the ideal of $A$ generated by all
products of the form $u_0(T_1^{a_1}\cdots T_r^{a_r})$, where
$(a_1,\dots,a_r)\in P_r$ is any sequence of exponents such
that $a_1+\cdots+a_r=\gamma$. The {\em angular blowing up}
of the ideal $T_\bullet A_0$ is the morphism
$$
Y:=\Proj\,(\sR\La A\Ra_{/\Q_+})\to X:=\Spec\,A
$$
(notation of definition \ref{def_Gamma-graded-algs}(v)).
Notice that $U_\gamma(\sR\La A\Ra_{/\Q_+})=Y$
for every $\gamma\in\Gamma\setminus\{0\}$.

\begin{theorem}\label{th_ceiling-cohomol}
In the situation of \eqref{subsec_angular-blow-up}, suppose
furthermore that the ring $A$ fulfills the condition of lemma
{\em\ref{lem_generalissimo}}. Then, for every $s\in\R_+$
we have :
\begin{enumerate}
\item
The natural maps
$T^{\lceil s\rceil}_\bullet A\to
T^{\lceil s\rceil}_\bullet\cdot H^0(Y,\cO_Y)\to
H^0(Y,T^{\lceil s\rceil}_\bullet\cO_Y)$
are isomorphisms.
\item
$H^p(Y,T^{\lceil s\rceil}_\bullet\cO_Y)=0$ for every $p>0$.
\item
The ring $\sR\La A\Ra$ fulfills the condition of lemma
{\em\ref{lem_generalissimo}}.
\item
For every $t\in\R_+$ the natural map
$$
T_\bullet^{\lceil s\rceil}A\otimes_AT_\bullet^{\lceil t\rceil}A\to
T_\bullet^{\lceil t+s\rceil}A
$$
is an isomorphism.
\end{enumerate}
\end{theorem}
\begin{proof} Notice first that
$T^{\lceil s\rceil}_\bullet\cdot H^0(Y,\cO_Y)$ and
$H^0(Y,T^{\lceil s\rceil}_\bullet\cO_Y)$ are $A$-submodules
of $ H^0(Y,\cO_Y)$; hence, the map
$T^{\lceil s\rceil}_\bullet\cdot H^0(Y,\cO_Y)\to
H^0(Y,T^{\lceil s\rceil}_\bullet\cO_Y)$ is injective, and in order
to show (i), it suffices to prove that the composition
$T^{\lceil s\rceil}_\bullet A\to H^0(Y,T^{\lceil s\rceil}_\bullet\cO_Y)$
is an isomorphism.

Now, for every $k\in\N$ we set
$\Gamma_k:=\{n/p^k~|~n\in\N\}\subset\Gamma$, and consider
the $\Gamma_k$-graded ring
$$
S_k:=\bigoplus_{\gamma\in\Gamma_k}I_k^{p^k\gamma}A
\qquad\text{and the $\Q_+$-graded ring}\qquad
S'_k:=(S_k)_{/\Q_+}
$$
where $I_k\subset R_{r,0}$ denotes the ideal generated by
$(T_1^{1/p^k},\dots,T_k^{1/p^k})$. The natural inclusion map
is a morphism of $\Q_+$-graded rings $S'_k\to S'_{k+1}$
for every $k\in\N$, and the colimit of the resulting
system $(S'_k~|~k\in\N)$ is $\sR(A)_{/\Q_+}$. Set as well
$$
Y_k:=\Proj\,S'_k
\qquad
\text{for every $k\in\N$}.
$$
There follows a system of affine morphisms of $X$-schemes
$(Y_k~|~k\in\N)$, and in view of
\cite[Ch.IV, Prop.8.2.3]{EGAIV-3} it is easily seen that
its limit is $Y$. For every $k\in\N$, let $\pi_k:Y\to Y_k$
be the natural projection. Notice that the direct summand
$S'_{k,\gamma}$ of $S'$ is an ideal of $A$, and
$$
\bigcup_{k\in\N}S'_{k,\gamma}=T_\bullet^{\La\gamma\Ra}A
\qquad
\text{for every $\gamma\in\Gamma$}.
$$
We deduce natural identifications as in remark
\ref{rem_blowing-up}(iv) for every $\gamma\in\Gamma$
and every $k\in\N$
$$
\cO_{Y_k}(\gamma)\isom S'_{k,\gamma}\cO_{Y_k}
\qquad
\colim_{i\in\N}\pi_i^{-1}\cO_{Y_i}(\gamma)\isom\cO_Y(\gamma)
\isom T_\bullet^{\La\gamma\Ra}\cO_Y
$$
whence, by proposition \ref{prop_dir-im-and-colim}(ii),
a natural isomorphism of $A$-modules :
\set\begin{equation}\label{eq_finish-by-virtue}
H^p(Y,T_\bullet^{\La\gamma\Ra}\cO_Y)\isom
\colim_{i\in\N}H^p(Y_i,\cO_{Y_i}(\gamma))
\qquad
\text{for every $p\in\N$ and every $\gamma\in\Gamma$}.
\end{equation}
By the same token, we have as well a natural isomorphism
\set\begin{equation}\label{eq_dir-lim-coh-angblowup}
H^p(Y,T_\bullet^{\lceil s\rceil}\cO_Y)\isom
\colim_{\gamma>s}H^p(Y,T_\bullet^{\La\gamma\Ra}\cO_Y)
\qquad
\text{for every $p\in\N$ and every $s\in\R_+$}.
\end{equation}

\begin{claim}\label{cl_universal-c}
There exists an integer $c>0$ such that for every
$k\in\N$ and every $\gamma\in\Gamma_k$, the inclusion
map $i:S'_{k,\gamma+c/p^k}\cO_{Y_k}\to S'_{k,\gamma}\cO_{Y_k}$
induces the zero map in cohomology :
$$
0=H^p(Y_k,i):H^p(Y_k,\cO_{Y_k}(\gamma+c/p^k))\to
H^p(Y_k,\cO_{Y_k}(\gamma))
\qquad
\text{for every $p>0$}
$$
and moreover the kernel and cokernel of the natural
morphism of inverse systems
$$
(S'_{k,n/p^k}\to H^0(Y_k,\cO_{Y_k}(n/p^k))~|~n\in\N)
$$
is uniformly essentially zero with step $\leq c$.
\end{claim}
\begin{pfclaim} Recall that $Y_k$ is the colimit of the
filtered system $(Y_k^{(\gamma)}~|~\gamma\in\Gamma_k)$ defined
as in \eqref{subsec_proj-with-Q-greaded}. By inspecting
the construction, it is easily seen that $Z_k:=Y_k^{(1/p^k)}$ is
the blowing up morphism of the ideal of $\cO_{\!X}$ generated
by the sequence $\bff^{(k)}:=(u_0(T_1^{1/p^k}),\dots,u_0(T_k^{1/p^k}))$.
Say now that $\gamma=a/p^k$; according to
\eqref{eq_divide-and-twist}, we have a natural identification
$$
\cO_{Y_k}(\gamma+d/p^k)\isom\Omega^{(1/p^k)}_*\cO_{Z_k}(a+d)
\qquad
\text{for every $d\in\N$}.
$$
On the other hand, lemma \ref{lem_generalissimo}(i)
implies that the ring $A$ satisfies condition
$\mathrm{(a)}^\mathrm{un}_{\bff^{(k)}}$ with step $\leq r$.
Then, by theorem \ref{th_too-difficult-for-me} and remark
\ref{rem_too-difficult-for-me}, there exists a constant
$c\in\N$ independent of $k$ and $a$, such that
the inclusion map $\cO_{Z_k}(a)\to\cO_{Z_k}(a+c)$ induces
the zero map in cohomology
$$
H^p(Z_k,\cO_{Z_k}(a+c))\to H^p(Z_k,\cO_{Z_k}(a))
\qquad
\text{for every $p>0$}
$$
and moreover, the kernel and cokernel of the morphism
of inverse systems
$$
(I^n_kA\to H^0(Z_k,\cO_{Z_k}(n))~|~n\in\N)
$$
are uniformly essentially zero with step $\leq c$.
Both assertions of the claim are an immediate consequence.
\end{pfclaim}

Now, in view of \eqref{eq_dir-lim-coh-angblowup}, in
order to prove (ii), it suffices to show that for every
$\gamma>s$ there exists $\gamma'\in\Q_+$ with
$s<\gamma'<\gamma$, and such that the natural map
$$
H^p(Y,T_\bullet^{\La\gamma\Ra}\cO_Y)\to
H^p(Y,T_\bullet^{\La\gamma'\Ra}\cO_Y)
$$
vanishes. To this aim, pick any $k\in\N$ such that
$c/p^k<\gamma-s$, where $c\in\N$ is as in claim
\ref{cl_universal-c}; by virtue of
\eqref{eq_finish-by-virtue}, we get the sought
vanishing with $\gamma':=\gamma-c/p^k$. A similar
argument proves also assertion (i) : details left
to the reader.

(iii): It suffices to check that for every $\gamma'<\gamma$
in $\N[1/p]$ the inclusion
$I_{r,0}^{\La\gamma\Ra}\subset I_{r,0}^{\La\gamma'\Ra}$ induces
the zero morphism
$$
\Tor_i^{R_{r,0}}(R_{r,0}/I_{r,0}^{\La\gamma\Ra},\sR\La A\Ra)\to
\Tor_i^{R_{r,0}}(R_{r,0}/I_{r,0}^{\La\gamma'\Ra},\sR\La A\Ra)
\qquad
\text{for every $i>0$}.
$$
However, for every $k\in\N$ set
$A_{r,k}:=\Z[T_1^{1/p^k},\dots,T_r^{1/p^k}]$ and let
$J_{r,k}\subset A_{r,k}$ be the ideal generated by
$(T_1^{1/p^k},\dots,T^{1/p^k}_r)$; we regard $A$ as an
$A_{r,k}$-algebra for every $k\in\N$, via restriction
of scalars along the inclusion map $A_{r,k}\to R_{r,0}$,
and we let $R_k(A)$ be the Rees $A$-algebra associated
with the $J_{r,k}A$-adic filtration of $A$. Now, say
that $\gamma=c/p^n$ and $\gamma'=c'/p^n$ for some
$c,c',n\in\N$; in view of claim \ref{cl_colim&Tors}(i)
it suffices to show that for every $i>0$ there exists
$k_0\in\N$ such that the inclusion
$J_{r,n+k}^{cp^k}\subset J_{r,n+k}^{c'p^k}$ induces the zero
morphism
$$
\Tor_i^{A_{r,n+k}}(A_{r,n+k}/J_{r,n+k}^{cp^k},R_{n+k}(A))\to
\Tor_i^{A_{r,n+k}}(A_{r,n+k}/J_{r,n+k}^{c'p^k},R_{n+k}(A))
$$
for every $k\geq k_0$. However, for every $k\in\N$ let
$\bff^{(k)}_\bullet:=(u_0(T_1^{1/p^k}),\dots,u_0(T^{1/p^k}_r))$,
and denote by $\bg^{(k)}_\bullet$ the image of $\bff^{(k)}_\bullet$
under the natural ring homomorphism $A\to R_k(A)$. By
lemma \ref{lem_generalissimo}(i) we know that $A$ fulfills
condition $\mathrm{(a)^{un}_{\bff^{(k)}_\bullet}}$ with step
$\leq r$; then corollary \ref{cor_a-unif-for-noether}(ii)
and remark \ref{rem_a-unif-for-noether} imply that $R_k(A)$
fulfills condition $\mathrm{(a)^{un}_{\bg^{(k)}_\bullet}}$ with a
step independent of $k$; explicitly, the latter means that
the inverse system
$(\Tor_i^{A_{r,k}}(A_{r,k}/J_{r,k}^t,R_k(A))~|~t\in\N)$ is
uniformly essentially zero for every $k\in\N$ and every
$i>0$, with a step independent of $k$. The assertion is
an immediate consequence.

(iv): It suffices to show that the natural map
\set\begin{equation}\label{eq_angle-suffice}
T_\bullet^{\lceil s\rceil}A\otimes_AT_\bullet^{\La t\Ra}A\to
T_\bullet^{\lceil t+s\rceil}A
\end{equation}
is an isomorphism for every $t\in\R_+$. However,
(iii) easily implies that the natural map
$$
T_\bullet^{\lceil s\rceil}R_{r,0}\otimes_{R_{r,0}}T_\bullet^{\La t\Ra}A
\to T_\bullet^{\lceil t+s\rceil}A
$$
is an isomorphism, and the latter factors through
\eqref{eq_angle-suffice} and a surjection
$T_\bullet^{\lceil s\rceil}R_{r,0}\otimes_{R_{r,0}}T_\bullet^{\La t\Ra}A
\to T_\bullet^{\lceil s\rceil}A\otimes_AT_\bullet^{\La t\Ra}A$,
whence the contention.
\end{proof}

\sset\subsubsection{}\label{subsec_situation-sit}
In the situation of \eqref{subsec_angular-blow-up},
suppose that $A$ fulfills the condition of lemma
\ref{lem_generalissimo}, and is endowed with a
$\Q_+$-grading which makes it an $A_0$-subalgebra
of $A_0[\Q_+]$ (where $A_0\subset A$ is the subring
of homogeneous elements of degree zero), and moreover
suppose that the image of $u_0$ lies in $A_0$. In this
case, the foregoing considerations apply to $A_0$ as well,
and we get a well defined $\Gamma$-graded subring
$\sR\La A_0\Ra$ of $\sR\La A\Ra$, as well the angular
blowing up morphism
$$
Y_0:=\Proj\,(\sR\La A_0\Ra_{/\Q_+})\to X_0:=\Spec\,A_0.
$$
Notice also that $\sR\La A\Ra$ carries a natural
$(\Gamma\times\Q_+)$-graded ring structure : namely,
$$
\gr_{(\gamma,t)}\sR\La A\Ra:=T^{\La\gamma\Ra}_\bullet A_t
\qquad
\text{for every $(\gamma,t)\in\Gamma\times\Q_+$}
$$
which is an ideal in $A_0$, for every such $(\gamma,t)$,
and we set
$$
\sR\La A_t\Ra:=\bigoplus_{\gamma\in\Gamma}T^{\La\gamma\Ra}_\bullet A_t=
A_t\sR\La A_0\Ra
\qquad
\text{for every $t\in\Q_+$}.
$$

\begin{corollary}\label{cor_ceiling-cohomol}
In the situation of \eqref{subsec_situation-sit}, for
every $s\in\R_+$ and every $t\in\Q_+$ we have :
\begin{enumerate}
\item
The natural maps $T^{\lceil s\rceil}_\bullet A_t\to
T^{\lceil s\rceil}_\bullet A_t\cdot H^0(Y_0,\cO_{Y_0})\to
H^0(Y_0,T^{\lceil s\rceil}_\bullet A_t\cO_{Y_0})$
are isomorphisms.
\item
$H^p(Y_0,T^{\lceil s\rceil}_\bullet A_t\cO_{Y_0})=0$
for every $p>0$.
\end{enumerate}
\end{corollary}
\begin{proof} Let $i:A_0\to A$ and
$\sR\La i\Ra:\sR\La A_0\Ra\to\sR\La A\Ra$ be the
inclusion maps; since $i(\sR\La A_0\Ra_+)$ generates the
ideal $\sR\La A\Ra_+$ of $\sR\La A\Ra$, then $i$ induces
a morphism of schemes $\Proj\,\sR\La i\Ra:Y\to Y_0$ (remark
\ref{rem_likewise}(ii)) fitting into the commutative diagram
$$
\xymatrix{ Y \ar[rr]^-{\Proj\,\sR\La i\Ra} \ar[d]_\pi & &
Y_0 \ar[d]^{\pi_0} \\
X \ar[rr]^-{\Spec\,i} & & X_0.
}$$
To ease notation, we let $\phi_X:=\Spec\,i$ and
$\phi_Y:=\Proj\,\sR\La i\Ra$. Notice that
$$
T^{\lceil s\rceil}_\bullet\sR\La A\Ra_{(f)}=
T^{\lceil s\rceil}_\bullet
\Bigl(\bigoplus_{t\in\Q_+}\sR\La A_t\Ra\Bigr)_{(f)}=
\bigoplus_{t\in\Q_+}T^{\lceil s\rceil}_\bullet A_t\sR\La A_0\Ra_{(f)}
$$
for every homogeneous element $f\in\sR\La A_0\Ra$,
which shows that
$$
\phi_{Y*}T^{\lceil s\rceil}_\bullet\cO_Y=
\bigoplus_{t\in\Q_+}T^{\lceil s\rceil}_\bullet A_t\cO_{Y_0}.
$$
Since $\phi_Y$ is an affine morphism, we deduce
$$
H^p(Y,T^{\lceil s\rceil}_\bullet\cO_Y)=
H^p(Y_0,\phi_{Y*}T^{\lceil s\rceil}_\bullet\cO_Y)=
\bigoplus_{t\in\Q_+}H^p(Y_0,T^{\lceil s\rceil}_\bullet A_t\cO_{Y_0})
$$
and combining with theorem \ref{th_ceiling-cohomol} we
get the corollary.
\end{proof}

\begin{example}\label{ex_special}
(i)\ \
According to proposition \ref{prop_back-to-toids}, theorem
\ref{th_ceiling-cohomol} applies especially to the case
where $A$ is perfectoid. In this case there exist unique
$\beta_1,\dots,\beta_r\in\bE:=\bE(A)$ such that
$u_0(T^{1/p^k}_i)=\bar u_A(\beta^{1/p^k}_i)$ for every
$i=1,\dots,r$ and every $k\in\N$, so we may regard as well
$\bE$ as an $R_{r,0}$-algebra via the unique ring homomorphism
$$
R_{r,0}\to\bE
\qquad :\qquad
T_i\mapsto\beta_i
\qquad
\text{for every $i=1,\dots,r$}.
$$

(ii)\ \
More generally, proposition \ref{prop_a-unif-for-perf} shows
that theorem \ref{th_ceiling-cohomol} applies to the case where
$A$ is the graded subring of a perfectoid ring $(A^\wedge,A)$
with $\Delta$-graded structure, where $\Delta$ is an integral
$p$-perfect monoid, provided the elements $\beta_1,\dots,\beta_r$
as in (i) are homogeneous elements of the $\Delta$-graded
subring of $\bE(A)$.

(iii)\ \
A ring $A$ as in (ii) is given by example \ref{ex_E-of-beta} :
we take a perfectoid ring $A_0$, a second finite family
$\beta'_\bullet:=(\beta'_1,\dots,\beta'_{r'})$ of elements
of $\bE_0:=\bE(A_0)$, we set $I:=\beta'_\bullet\bE_0$, and
let $A_t$ be the topological closure in $A_0$ of $I^{\La t\Ra}A_0$
for every $t\in\Delta:=\N[1/p]$. Then, clearly corollary
\ref{cor_ceiling-cohomol} shall apply to the $\Q_+$-graded
ring $A_{/\Q_+}$.

(iv)\ \
Let $A_0$, $\bE_0$, $\beta'_\bullet$ and $A$ be as in (iii), and
suppose that also the sequence $\beta_\bullet$ lies in $\bE(A_0)$.
Suppose moreover that $I$ is an ideal of adic definition for
$\bE_0$. In this case, the ideal $I^{\La t\Ra}A_0$ is also
open in $A_0$ (corollary \ref{cor_taut-two}(ii)), and therefore
it coincides with $A_t$ for every $t\in\N[1/p]$. Especially,
$\gr_{(\gamma,t)}\sR\La A\Ra=T^{\La\gamma\Ra}_\bullet I^{\La t\Ra}A_0$
for every $(\gamma,t)\in\Gamma\times\N[1/p]$.

(v)\ \
In the situation of (iv), suppose furthermore that the
ideal $T_\bullet\bE_0\subset\bE_0$ generated by $\beta_\bullet$
is open. Notice that $T^{\lceil 0\rceil}_\bullet\bE_0$ is the
radical of $T_\bullet\bE_0$ (details left to the reader).
It follows that for every $s\in\N[1/p]\setminus\{0\}$
there exists $s'>0$ such that
$I^{\La s\Ra}\bE_0\subset T_\bullet^{\La s'\Ra}\bE_0$, whence
$I^{\La s\Ra}A_0\subset T_\bullet^{\La s'\Ra}A_0$ (corollary
\ref{cor_bingo}(ii)). Therefore :
$$
I^{\La t+s\Ra}A_0\subset
I^{\La t\Ra}T_\bullet^{\lceil 0\rceil}A_0\subset I^{\La t\Ra}A_0
\qquad
\text{for every $t,s\in\N[1/p]$ with $s>0$}.
$$
\end{example}

\begin{corollary}\label{cor_angular-blow-up}
In the situation of \eqref{subsec_angular-blow-up},
suppose that $A$ is perfectoid, and let $I$ be any
finitely generated ideal of adic definition of\/
$\bE:=\bE(A)$. With the notation of example
{\em\ref{ex_special}(v)}, for every $t\in\R_+$ we have :
\begin{enumerate}
\item
The natural maps
$I^{\lceil t\rceil}A\to I^{\lceil t\rceil}\cdot H^0(Y,\cO_Y)
\to H^0(Y,I^{\lceil t\rceil}\cO_Y)$ are isomorphisms.
\item
$H^p(Y,I^{\lceil t\rceil}\cO_Y)=0$ for every $p>0$.
\end{enumerate}
\end{corollary}
\begin{proof} With the notation of example
\ref{ex_special}(v), notice that
$$
H^p(Y,I^{\lceil t\rceil}T_\bullet^{\lceil 0\rceil}\cO_Y)=
\colim_{s>t}H^p(Y,I^{\La s\Ra}T_\bullet^{\lceil 0\rceil}\cO_Y).
$$
On the other hand, example \ref{ex_special}(v) implies that
$I^{\lceil t\rceil}A=I^{\lceil t\rceil}T_\bullet^{\lceil 0\rceil}A$
for every $t\in\R_+$. To conclude, it suffices to invoke
corollary \ref{cor_ceiling-cohomol} as in example
\ref{ex_special}(iii,iv) : {\em i.e.} with $A_0$ and $Y_0$
replaced by the current $A$ and respectively $Y$, and with
$A$ replaced by the $\Q_+$-graded ring $A'_{/\Q_+}$,
where $A'_s:=I^{\La s\Ra}A$ for every $s\in\N[1/p]$.
\end{proof}

\begin{definition}\label{def_formal-perfectoid}
We say that a topological ring $(A,\cT)$ is a
{\em formal P-ring} (resp. a {\em formal perfectoid ring})
if it is adic with a finitely generated ideal of adic
definition, and the separated completion of $(A,\cT)$
is a P-ring (resp. a perfectoid ring).
\end{definition}

\begin{remark} Let $(A,\cT)$ be any topological ring.

(i)\ \
Suppose that $(A,\cT)$ is a P-ring. Then, any ideal
of definition of the P-ring $A$ is also an ideal of adic
definition of the adic ring $A$. Remark \ref{rem_why-only-now}
implies that the converse holds if and only if $A$ is an
$\F_p$-algebra. Indeed, if $I$ is any ideal of adic
definition of $A$, then any power $I^n$ is also an
ideal of adic definition, but if $p$ does not vanish
in $A$, then one such power will not contain $p$,
whence the contention.

(ii)\ \
Let $\phi:(A,\cT)\to(A',\cT')$ be an adic morphism of
topological rings, and suppose that $\phi$ is a weakly
\'etale ring homomorphism. Then \cite[Lemma 3.1.2(i)]{Ga-Ra}
implies that $\phi$ is adically weakly \'etale (see definition
\ref{def_c-adic-and-adic}(iii)).
\end{remark}

\begin{theorem}\label{th_formal-perf}
Let $f:A\to B$ be a c-adically weakly \'etale morphism
of topological rings.
\begin{enumerate}
\item
If $A$ and $B$ are f-adic and $A$ is a formal P-ring
(resp. a formal perfectoid ring), the same holds for $B$.
\item
If $f$ is c-adically faithfully flat, the following holds :
\begin{enumerate}
\item
If $B$ is f-adic and $A$ and $B$ are complete and
separated, $A$ is f-adic and $f$ is adic.
\item
If $B$ is a formal P-ring (resp. a formal perfectoid ring),
the same holds for $A$.
\end{enumerate}
\end{enumerate}
\end{theorem}
\begin{proof}(i): In light of lemma \ref{lem_still-c-adic}(iii,iv),
we may assume that $A$ and $B$ are complete and separated,
in which case $A$ is a P-ring (resp. a perfectoid ring),
$f$ is adically weakly \'etale (lemma \ref{lem_f-adics}(i.c)),
and we have to show that $B$ is a P-ring (resp. is perfectoid).

Suppose first that $A$ is a P-ring, and let $I$ be any
ideal of definition of $A$. It suffices to check that
$IB$ is an ideal of definition for $B$. However, the
natural map $B/I^2B\to B^\wedge/I^2B^\wedge$ is an isomorphism,
by remark \ref{rem_completion-of-topring}(ii), so we come
down to checking that the Frobenius endomorphism of $B/I^2B$
is surjective. However, the induced ring homomorphism
$A/I^2\to B/I^2B$ is weakly \'etale by assumption; since
the Frobenius endomorphism of $A/I^2$ is surjective by
assumption, the assertion then follows easily from
\cite[Th.3.5.13(ii)]{Ga-Ra}.

Next, suppose that $A$ is perfectoid. Due to remark
\ref{rem_p-can-lie-deep}, we may assume that $p\in I^t$,
where $t\in\N$ is defined as in
\eqref{subsec_setup-critperf}, in terms of the
length of a system of generators for $I$. Set
$J:=I^{(p)}$; taking into account (i) and theorem
\ref{th_criterium-perfect}, we are then reduced to
checking that the morphism of graded rings
$$
\Phi_{IB}:\gr^\bullet_{IB}B\to
\gr^\bullet_{JB}B
$$
is an isomorphism. However, since $f\otimes_AA/I^n$
is flat for every $n\in\N$, we get induced graded
ring isomorphisms
$$
\phi_I:B\otimes_A\gr^\bullet_IA\isom\gr^\bullet_{IB}B
\qquad
\phi_J:B\otimes_A\gr^\bullet_JA\isom\gr^\bullet_{JB}B
$$
and $\phi_I$ identifies the map
$$
\gr^\bullet_If:\gr^\bullet_IA\to\gr^\bullet_{IB}B
$$
with the base change $f\otimes_A\gr^\bullet_IA$.
Especially, since $f\otimes_AA/I$ is weakly \'etale,
the same holds for $\gr^\bullet_If$
(\cite[Lemma 3.1.2(i)]{Ga-Ra}). Moreover, denote by
$(\gr^\bullet_IA)_{(\Phi)}$ the $\gr^\bullet_IA$-algebra
whose underlying ring is $\gr^\bullet_IA$ and whose
structure map is the Frobenius endomorphism
$\Phi_{\gr^\bullet_IA}$, and define likewise
$(\gr^\bullet_{IB}B)_{(\Phi)}$; by
\cite[Th.3.5.13(ii)]{Ga-Ra}, the maps $\phi_J$ and
$\gr^\bullet_Jf$ induce an isomorphism
$$
\psi_J:B\otimes_A(\gr^\bullet_JA)_{(\Phi)}\isom
\gr^\bullet_{JB}B\otimes_{\gr^\bullet_JA}
(\gr^\bullet_JA)_{(\Phi)}
\isom(\gr^\bullet_{JB}B)_{(\Phi)}
\qquad
b\otimes a\mapsto b^p\cdot\gr_If(a).
$$
Summing up, we get a commutative diagram of ring
homomorphisms
\set\begin{equation}\label{eq_long-day}
{\diagram
B\otimes_A\gr^\bullet_IA \ar[rr]^-{B\otimes_A\Phi_I}
\ar[d]_{\phi_I} & &
B\otimes_A(\gr^\bullet_JA)_{(\Phi)} \ar[d]^{\psi_J} \\
\gr^\bullet_{IB}B \ar[rr]^-{\Phi_{IB}} & &
(\gr^\bullet_{JB}B)_{(\Phi)}
\enddiagram}
\end{equation}
whose vertical arrows are isomorphisms.
Lastly, since $A$ is a formal perfectoid ring, $\Phi_I$
is an isomorphism as well (theorem \ref{th_criterium-perfect}),
whence the assertion.

(ii.a): Let $J$ be any finitely generated ideal of adic
definition of $B$, and $(I_\lambda~|~\lambda\in\Lambda)$ a
cofiltered system of ideals that defines the linear topology
of $A$. By assumption, $(I_\lambda B)^c$ is open in $B$,
hence $I_\lambda B=(I_\lambda B)^c$ for every
$\lambda\in\Lambda$ (lemma \ref{lem_5.3.8}(ii.b));
then the induced map
$$
f_\lambda:A_\lambda:=A/I_\lambda\to B_\lambda:=B/I_\lambda B
$$
is weakly \'etale and faithfully flat for every
$\lambda\in\Lambda$. Also by assumption, we may find
$\mu\in\Lambda$ such that $I_\mu B\subset J^2$.
Hence, $JB_\mu$ is contained in the nilradical
$\nil(B_\mu)$ of $B_\mu$, and it follows easily
from claim \ref{cl_generically-true}(i) that
$$
\nil(B_\mu)=\nil(A_\mu)\cdot B_\mu.
$$
Therefore, we may find a finitely generated ideal
$K'\subset\nil(A_\mu)$ such that $JB_\mu\subset K'B_\mu$.
Pick a finitely generated ideal $K\subset A$ such that
$KA_\mu=K'$.
It follows easily that $J\subset KB+J^2$. However,
$B$ is complete and $J$ is topologically nilpotent
in $B$, so $J$ lies in the Jacobson radical of $B$,
and consequently $J\subset KB$, by Nakayama's lemma.
On the other hand, by construction, we may find $n\in\N$
such that $K'{}^n=0$, whence $K^n\subset I_\mu$, and
therefore $K^nB\subset J^2$. Summing up, we conclude
that $KB$ is an ideal of adic definition of $B$.
For every $n\in\N$, pick $\nu(n)\in\Lambda$ such that
$I_{\nu(n)}B\subset J^n$; then
$I_{\nu(n)}B_\lambda\subset K^nB_\lambda$ for every
$\lambda\in\Lambda$, so that
$I_{\nu(n)}A_\lambda\subset K^nA_\lambda$, as $f_\lambda$
is faithfully flat. We deduce that the topological
closure $(K^n)^c$ of $K^n$ contains $I_{\nu(n)}$, and
especially, it is an open ideal of $A$, for every
$n\in\N$. On the other hand, for every $\lambda\in\Lambda$
we may find $k(\lambda)\in\N$ such that
$J^{2k(\lambda)}\subset I_\lambda B$, so that
$K^{n\cdot k(\lambda)}B_\nu\subset I_\lambda B_\nu$
for every $\nu\geq\lambda$. It follows that
$K^{n\cdot k(\lambda)}A_\nu\subset I_\lambda A_\nu$,
again since $f_\nu$ is faithfully flat; then, since
$I_\lambda$ is closed in $A$, we deduce that
$K^{n\cdot k(\lambda)}\subset I_\lambda$, for every
$\lambda\in\Lambda$. Summing up, we see that
the topology of $A$ is c-adic, hence f-adic, by
lemma \ref{lem_5.3.8}(ii.a). Then $f$ is adic,
by lemma \ref{lem_f-adics}(i.c).

(ii.b): In light of lemma \ref{lem_still-c-adic}(iii), we
may replace $f$ by $f^\wedge$, and assume from start that
$A$ and $B$ are complete and $B$ is a P-ring (resp. is
perfectoid), and we need to show that the same holds for
$A$. However, (ii.a) already says that $A$ is f-adic, and
$f$ is adic. Suppose first that $B$ is a P-ring, and pick
any finitely generated ideal of adic definition for $A$;
without loss of generality we may assume that $J$ is an
ideal of definition for the P-ring $B$, and that
$IB\subset J$. Set $A_n:=A/I^{n+1}$ and $B_n:=B/I^{n+1}B$
for every $n\in\N$, and notice that $JB_0$ is nilpotent
ideal of $B_0$; arguing as in the proof of (ii.a), we
deduce that there exists a finitely generated ideal
\set\begin{equation}\label{eq_put-it-here}
K_0\subset\nil(A_0)
\end{equation}
such that $JB_0\subset K_0B_0$. Then, the preimage $K$
of $K_0$ in $A$ is a finitely generated ideal such that
$JB\subset KB$; hence $p^2\in K^2B_n$ for every $n\in\N$.
Since the induced map $A_n\to B_n$ is faithfully flat,
it follows that $p\in K^2A_n$ for every $n\in\N$, so
$p$ lies in the topological closure of $K^2$ in $A$;
but $(K^2)^c=K^2$, since $K$ is open in $A$ by
construction. Lastly, \eqref{eq_put-it-here} implies
that $K^n\subset I$ for some sufficiently large $n\in\N$;
thus, $K$ is also an ideal of adic definition of $A$,
and we are reduced to checking that the Frobenius
endomorphism $\Phi_{A/J^2}$ of $A/J^2$ is surjective.
However, by lemma \ref{lem_perfectoid}(iv) the
Frobenius endomorphism $\Phi_{B/J^2B}$ is surjective;
by \cite[Th.3.5.13(ii)]{Ga-Ra}, it follows that
$\Phi_{A/J^2}\otimes_{A/J^2}B/J^2B$ is also surjective.
But the induced map $A/J^2\to B/J^2B$ is faithfully
flat, whence the contention.

Lastly, suppose that $B$ is perfectoid; by the foregoing,
we know already that $A$ is a P-ring, and by remark
\ref{rem_p-can-lie-deep} we may find an ideal of
definition $I$ of $A$ such that $p\in I^t$, where
$t\in\N$ is defined as in \eqref{subsec_setup-critperf}.
We set $J:=I^{(p)}$; arguing as in the proof of (i), we
get a commutative diagram \eqref{eq_long-day}, whose
vertical arrows are isomorphisms. Now, notice that,
since $f$ is adic, $IB$ is an ideal of definition of
$B$, so the bottom horizontal arrow of \eqref{eq_long-day}
is an isomorphism; hence the same holds for
$B\otimes_A\Phi_I=(B/IB)\otimes_{A/I}\Phi_I$. But by
assumption, $f$ induces a faithfully flat map
$A/I\to B/IB$, so $\Phi_I$ is an isomorphism, and
therefore $A$ is perfectoid, by theorem
\ref{th_criterium-perfect}.
\end{proof}

\sset\subsubsection{}\label{subsec_Letta-caduto}
Let $\underline U:=(U,\cT_A,A^+_U)$ be any perfectoid
quasi-affinoid scheme, and set
$\underline U{}_\bE:=\bE(\underline U)$ (notation of
\eqref{subsec_upgrade-E}), so that
$\underline U{}_\bE=(\bE(U),\cT_\bE,\bE^+_U)$, where
$\bE^+_U:=\bE(A^+_U)$. Let as well
$$
A_U:=\cO_{\!U}(U)
\qquad
\bE_U:=\cO_{\bE(U)}(\bE(U))
$$
and pick a perfectoid subring of definition $A$ of $A_U$
{\em such that $A^+_U$ is the integral closure of $A\cap A^+_U$
in $A_U$ and such that the induced map $U\to\Spec\,A$ is an
open immersion} (by theorem \ref{th_int-subrings-perfectoid}(i)
and lemma \ref{lem_dots-on-is}(i), the ring $A^\circ_U$ is one
such subring of definition), and recall that $\bE:=\bE(A)$ is
a ring of definition of $\bE_U$. We denote as usual by
$\phi_U^\flat:\bE_U\to A_U$ the continuous map of monoids
provided by proposition \ref{prop_new-formula}(i), and
recall that $\phi^\flat_U$ restricts to a continuous map
$\bar u_A:\bE\to A$. Now, let $e_0,\dots,e_n$ be a finite
system of elements of $\bE_U$ that generates an open ideal;
according to claim \ref{cl_extend-cor-taut-two}, the system
$(a_i:=\phi^\flat_U(e_i)~|~i=0,\dots,n)$ generates an open
ideal of $A_U$, so we may consider the rational subset
$$
R:=R_{A_U}\Bigl(\frac{a_1}{a_0},\dots,\frac{a_n}{a_0}\Bigr)\cap
\Spa\,\underline U
$$
and the corresponding topological rings $\cO^\wedge_{\Spa\,U}(R)$
and $\cO^{\wedge+}_{\Spa\,U}(R)$, defined as in
\eqref{subsec_rational-site}. Recall the construction of
$\cO^\wedge_{\Spa\,U}(R)$ : first, there is a natural f-adic
topology on the localization $A_R:=A_U[1/a_0]$ such that
the localization map $A_U\to A_R$ is f-adic, and the subring
$$
B_R:=A[a_i/a_0~|~i=1,\dots,n]\subset A_R
$$
is a ring of definition of $A_R$; hence $A_R$ induces on
$B_R$ the unique linear topology such that the natural map
$A\to B_R$ is adic. We let $A^+_R$ be the integral closure
of $A^+_U[a_i/a_0~|~i=1,\dots,n]$ in $A_R$, set
$U_R:=U\cap\Spec\,A_R$, and define
$$
\underline A{}_R:=(A_R,A^+_R,U_R)
\qquad\text{and}\qquad
\underline U_R^\wedge:=(\sSpec\,\underline A{}_R)^\wedge.
$$
Then (up to natural isomorphism) we have
$$
\sGamma(\underline U_R^\wedge)=
(\cO^\wedge_{\Spa\,U}(R),\cO^{\wedge+}_{\Spa\,U}(R),U^\wedge_R)
\qquad\text{where}\qquad
U^\wedge_R:=U_R\times_{\Spec\,A_R}\Spec\,\cO^\wedge_{\Spa\,U}(R).
$$
Notice that $A_R$ contains also the subrings
$$
\begin{aligned}
C_R:=&\,
A\Bigl[\phi^\flat_U(e_i^{1/p^k})/\phi^\flat_U(e_0^{1/p^k})~|~
(i,k)\in\{1,\dots,n\}\times\N\Bigr] \\
D_R:=&\,
A^+_U\Bigl[\phi^\flat_U(e_i^{1/p^k})/\phi^\flat_U(e_0^{1/p^k})~|~
(i,k)\in\{1,\dots,n\}\times\N\Bigr]
\end{aligned}
$$
and endow $C_R$ and $D_R$ with the topologies induced by
the inclusion into $A_R$. Let also $C_R^+:=C_R\cap A^+_R$.
The first observation is the following :

\begin{proposition}\label{prop_Letta-non-caduto}
With the notation of \eqref{subsec_Letta-caduto}, the
following holds :
\begin{enumerate}
\item
The natural morphism
$\underline A{}_R^\wedge\to\sGamma(\underline U_R^\wedge)$
is an isomorphism.
\item
The datum $\underline C{}_R:=(C_R,C^+_R,U_R)$ is a
quasi-affinoid ring, and the natural morphism
$\sSpec\,\underline A{}_R\to\sSpec\,\underline C{}_R$
is an isomorphism.
\item
$\underline C{}^\wedge_R$ is a perfectoid quasi-affinoid
ring, and $\underline U_R^\wedge$ is a perfectoid
quasi-affinoid scheme.
\item
Suppose moreover that 
$\beta_{\underline U}(U_R)\subset\Omega_{\underline U}$
(notation of definition {\em\ref{def_spread}}).
Then the datum $\underline D{}_R:=(D_R,D_R,U_R)$ is a
quasi-affinoid ring, and $\underline D{}_R^\wedge$ is
perfectoid.
\end{enumerate}
\end{proposition}
\begin{proof}(ii): We consider the natural morphisms of
schemes
$$
U_R\xrightarrow{\ f\ }\Spec\,A_R\xrightarrow{\ g\ }\Spec\,C_R
\xrightarrow{\ h\ }\Spec\,A
$$
whose composition is also the composition of the open
immersion $U_R\to U$ with the natural morphism
$i:U\to\Spec\,A$. By assumption, $i$ is an open immersion,
and its image contains the analytic locus of
$\Spec\,A$ (lemma \ref{lem_deja-vu}(iii)). It
follows that $h\circ g\circ f$ is an open immersion;
moreover, by construction, both $f$ and $g$ have
schematically dense images. From claim \ref{cl_reinstated},
we deduce that $g\circ f$ is an open immersion, and
$f(U_R)=(h\circ g)^{-1}i(U_R)=(h\circ g)^{-1}i(U)$.
Now, $U$ contains the analytic locus of $\Spec\,A$,
and the natural map $A\to A_R$ is f-adic; taking into
account lemma \ref{lem_deja-vu}(iv), it follows that
$f(U_R)$ contains the analytic locus of $\Spec\,A_R$.
Lastly, $g$ identifies the analytic locus of $\Spec\,A_R$
with that of $\Spec\,C_R$ (lemma \ref{lem_deja-vu}(iii)),
so $g\circ f(U_R)$ contains the analytic locus of $\Spec\,C_R$.
Summing up, this shows that $\underline C{}_R$ is a
quasi-affinoid ring. Then, a simple inspection shows
that the inclusion map $C_R\to A_R$ induces an
isomorphism
$\sSpec\,\underline A{}_R\to\sSpec\,\underline C{}_R$.

(i): To begin with, we remark :

\begin{claim}\label{cl_add-roots}
$C_R$ and $D_R$ are subrings of definition of $A_R$.
\end{claim}
\begin{pfclaim} Clearly $C_R$ and $D_R$ are open in $A_R$,
and $D_R\subset C_R$; by virtue of proposition
\ref{prop_f-adics}(ii), we have then only to check that
$C_R$ is bounded in $A_R$. To this aim, since $B_R$ is
open and bounded in $A_R$, it suffices to show that
there exists an open ideal $J$ of $B_R$ such that
$J\cdot C_R\subset B_R$. However, pick any finite system
$f_1,\dots,f_r$ of generators for an ideal of definition
of $\bE$; since the system $e_0,\dots,e_n$ generates an open
ideal of $\bE_U$, we may find an integer $k\in\N$ large enough
so that the $f^r_ie_j\in\bE$ for every $i=1,\dots,r$ and every
$j=0,\dots,n$, and then it is clear that the system
$(f^r_ie_j~|~i\leq r,\ j\leq n)$ generates an open ideal of
$\bE$. By claim \ref{cl_extend-cor-taut-two}, it follows that
the system $(\phi^\flat_U(f^r_ie_j)~|~i\leq r,\ j\leq n)$ generates
an open ideal $J_0$ of $A$, and therefore $J:=J_0B_R$ is
open in $B_R$. To conclude, it suffices to show that
$\phi^\flat_U(f^r_ie_je^\nu_se^{-\nu}_0)\in B_R$ for every
$i=1,\dots,r$ every $j,s=0,\dots,n$, and every $\nu\in\N[1/p]$.
Moreover, we may easily reduce to the case where $0<\nu<1$.
Then we see that
$$
f^r_ie_je^\nu_se^{-\nu}_0=
\frac{e_j}{e_0}\cdot(f^{r\nu}_ie_s^\nu)\cdot(f^{r(1-\nu)}_ie_0^{1-\nu})
$$
where both $f^{r\nu}_ie_s^\nu$ and $f^{r(1-\nu)}_ie_0^{1-\nu}$ lie
in $\bE$, since the latter is a perfect $\F_p$-algebra. The
assertion follows.
\end{pfclaim}

Pick a finite system $\beta_\bullet:=(\beta_1,\dots,\beta_s)$
of elements of $\bE$ that generates an ideal of adic definition,
and set $f_i:=\bar u_A(\beta_i)$ for $i=1,\dots,s$, so that the
system $\bff_\bullet:=(f_1,\dots,f_s)$ generates an ideal
of adic definition of $A$. By claim \ref{cl_add-roots},
it then follows that the system $\bff_\bullet$ also generates
an ideal of adic definition $J$ for $C_R$. Denote by
$C_R^\wedge$ the separated completion of $C_R$; we notice :

\begin{claim}\label{cl_think-of-llamas}
The ring $C_R$ fulfills condition $\mathrm{(a)^{un}_\bff}$
of \eqref{subsec_badabum}, and $C^\wedge_R$ is perfectoid.
\end{claim}
\begin{pfclaim} Set $\Delta:=\N[1/p]$; recall that example
\ref{ex_E-of-beta} attaches to  the sequence
$e_\bullet:=(e_0,\dots,e_n)$ of elements of $\bE_U$, two
perfectoid rings with $\Delta$-graded structure
$$
(\cE(e_\bullet)^\wedge,\cE(e_\bullet))
\qquad\text{and}\qquad
(\cA(e_\bullet)^\wedge,\cA(e_\bullet)).
$$
Moreover, $\gr_0\cE(e_\bullet)=\bE$, $\gr_0\cA(e_\bullet)=A$,
and we have a natural isomorphism
$$
\bE(\cA(e_\bullet)^\wedge,\cA(e_\bullet))\isom
(\cE(e_\bullet)^\wedge,\cE(e_\bullet)).
$$
Furthermore, since $e_\bullet$ generates an open ideal of
$\bE_U$, the $\bE$-module $\gr_\delta\cE(e_\bullet)$ (resp.
the $A$-module $\gr_\delta\cA(e_\bullet)$) is the submodule
of $\bE_U$ (resp. of $A_U$) generated by the elements of
the form $e_0^{\nu_0}\cdots e_n^{\nu_n}$ (resp.
$\phi^\flat_U(e_0^{\nu_0}\cdots e_n^{\nu_n})$), where
$(\nu_0,\dots,\nu_n)$ ranges over all elements of
$\Delta^{\oplus n+1}$ such that $\nu_0+\cdots+\nu_n=\delta$
(see remark \ref{rem_beta-taut}(viii)). Especially, we
may regard $\beta_\bullet$ as a sequence of elements of
$\gr_0\cE(e_\bullet)$, whose image under the map
$\bar u_{\cA(e_\bullet)}:\cE(e_\bullet)\to\cA(e_\bullet)$ is
(naturally identified with) the sequence $\bff_\bullet$.
Next, let us endow $\cA':=\cA(e_\bullet)[a_0^{-1}]$ with
its $J\cA'$-adic topology; since we are inverting an
element $a_0\in\gr_1\cA(e_\bullet)$,  the $A$-algebra
$\cA'$ is also $\Delta$-graded, and a simple inspection
yields a natural isomorphism of $C_R\isom\gr_0\cA'$
of topological rings, which identifies the sequences
$\bff_\bullet$ in these two rings. Lastly, by virtue of
proposition \ref{prop_a-unif-for-perf}(ii), we deduce that
$\cA(e_\bullet)$ satisfies condition $\mathrm{(a)^{un}_\bff}$
of \eqref{subsec_badabum}, and then the same holds for
$\gr_0\cA'$, whence the first assertion.

In order to show that $C^\wedge_R$ is perfectoid, notice that
$J^n\cA'=\bigoplus_{\delta\in\Delta}J^n\gr_\delta\cA'$ for every
$n\in\N$. It follows easily that the maximal separated
quotient $\cA''$ of $\cA'$ is also a $\Delta$-graded ring
(so that the projection $\cA'\to\cA''$ is a homomorphism of
$\Delta$-graded rings), and $\gr_0\cA''$ is the maximal
separated quotient of $\gr_0\cA'$ : details left to the
reader. Hence $(\cA'',\cA'')$ is a topological ring with
$\Delta$-pre-graded structure, and therefore
$(\cC^\wedge,\cC):=(\cA'',\cA'')^\wedge$ is a topological
ring with $\Delta$-graded structure (proposition
\ref{prop_Cauchy}(iii)), and the induced map
$C_R^\wedge\to\gr_0\cC$ is an isomorphism of topological
rings. On the other hand, since $\cA(e_\bullet)^\wedge$ is
perfectoid, and since the localization map
$\cA(e_\bullet)\to\cA'$ is an adic ring homomorphism,
theorem \ref{th_formal-perf}(i) implies that $\cC^\wedge$
is perfectoid as well. To conclude, it suffices now to
invoke proposition \ref{prop_graded-perfectoid}(ii).
\end{pfclaim}

As the morphism $i:U\to\Spec\,A$ is an open immersion whose
image contains the analytic locus of $\Spec\,A$, we may find
an open ideal $I\subset A$ such that
$\Spec\,A\setminus U=\Spec\,A/I$, and we pick a finite system
$\bg_\bullet:=(g_1,\dots,g_r)$ of generators of $I$. On the other
hand, $i$ factors through the natural open immersion
$j:U\to\Spec\,A_U$ and the morphism $\phi:\Spec\,A_U\to\Spec\,A$
induced by the inclusion map $A\to A_U$; by claim
\ref{cl_reinstated}(ii), it follows that $\phi^{-1}i(U)=j(U)$,
whence $\Spec\,A_U/IA_U=\Spec\,A_U\setminus j(U)$, and therefore
$$
Z:=\Spec\,A_R/IA_R=\Spec\,A_R\setminus U_R.
$$
In view of claim \ref{cl_add-roots}, remark
\ref{rem_depth-and-completion}(ii) and proposition
\ref{prop_complete-f-adic}(iii), we are then reduced to
showing that
$$
\depth_ZA_R\otimes_{C_R}C^\wedge_R>1.
$$
However, set $Q^\bullet:=\Cone(C_R[0]\to C^\wedge_R[0])$; clearly
we have $\depth_{IA_U}A_U>1$, whence $\depth_ZA_R>1$, and thus
it suffices to check that
$$
R\Gamma_{\!Z}(A_R\otimes_{C_R}Q^\bullet)=0.
$$
In light of proposition \ref{prop_depth-Kosz}(i), we are
then reduced to showing that the Koszul complex
$\bK^\bullet(\bg_\bullet,A_R\otimes_{C_R}Q^\bullet)$ is acyclic.
We remark :

\begin{claim}\label{cl_same-as-derotimes}
The natural map
$A_R\derotimes_{C_R}Q^\bullet\to A_R\otimes_{C_R}Q^\bullet$ is
an isomorphism in $\sD(C_R\Mod)$.
\end{claim}
\begin{pfclaim} We have a natural morphism of distinguished
triangles
$$
\xymatrix{ A_R[0] \ar[r] \ddouble &
A_R[0]\derotimes_{C_R}C_R^\wedge[0] \ar[r] \ar[d] &
A_R\derotimes_{C_R}Q^\bullet \ar[r] \ar[d] & A_R[1] \ddouble \\
A_R[0] \ar[r] & A_R\otimes_{C_R}C_R^\wedge[0]
\ar[r] & A_R[0]\otimes_{C_R}Q^\bullet \ar[r] & A_R[1]
}$$
in light of which, we are reduced to checking that
$T_i:=\Tor^{C_R}_i(A_R,C_R^\wedge)=0$ for every $i>0$. However,
clearly the natural map $T_i\to\Tor^{C_R}_i(A_R/C_R,C_R^\wedge)$
is an isomorphism for $i>1$, and is injective for $i=1$, and
moreover $A_R/C_R=\bigcup_{n\in\N}\Ann_{A_R/C_R}(J^n)$. Hence, the
assertion follows from claims \ref{cl_more-will-vanish}
and \ref{cl_think-of-llamas}.
\end{pfclaim}

Taking into account lemma \ref{lem_koszul-vanish}(v) and
claim \ref{cl_same-as-derotimes}, we are further reduced
to showing that
$\bK_\bullet(\bg_\bullet)\derotimes_{C_R}A_R\derotimes_{C_R}Q^\bullet$
is acyclic, and to this aim it suffices to check that the
same holds for the complex
$\bK_\bullet(\bg_\bullet)\otimes_{C_R}Q^\bullet$. The latter
follows easily from lemma \ref{lem_koszul-vanish}(ii),
claim \ref{cl_think-of-llamas} and corollary
\ref{cor_cond-a-and-compl}(i).

(iii): From claim \ref{cl_think-of-llamas} we know already
that $\underline C{}^\wedge_R$ is perfectoid, and from (i) and
(ii) we get isomorphisms of quasi-affinoid schemes
$$
\underline U{}^\wedge_R\isom\sSpec\,(\underline A{}^\wedge_R)\isom
\sSpec\,(\underline C{}_R^\wedge)
$$
whence the contention.

(iv): We consider the natural morphisms of schemes :
$$
U_R\xrightarrow{\ f\ }\Spec\,A_R\xrightarrow{\ g\ }
\Spec\,D_R\to\Spec\,A^+_U
$$
whose composition is the same as the restriction of
$\beta_{\underline U}$ to $U_R$, hence it is an open
immersion. Notice also that $\beta_{\underline U}$ is in turn
the composition of the open immersion $i:U\to\Spec\,A$
and the natural morphism $\Spec\,A\to\Spec\,A^+_U$;
the latter identifies the analytic locus of $\Spec\,A$
with that of $\Spec\,A^+_U$ (lemma \ref{lem_deja-vu}(iii)),
hence the image of $\beta_{\underline U}$ contains the
analytic locus of $\Spec\,A^+_U$. We may then argue
as in the proof of (ii) to conclude that $g\circ f$
is an open immersion and $g\circ f(U_R)$ contains the
analytic locus of $\Spec\,D_R$. Together with claim
\ref{cl_add-roots}, this already shows that
$\underline D{}_R$ is a quasi-affinoid ring. Moreover,
we also know that $A^+_U$ is a perfectoid ring (theorem
\ref{th_int-subrings-perfectoid}(iii)), and our
assumptions imply that $\underline A':=(A^+_U,A^+_U,U)$
is another perfectoid quasi-affinoid ring with
$\sSpec\,\underline A'=\underline U$. We may then
argue as in the proof of claim \ref{cl_think-of-llamas}
to see that the separated completion $D^\wedge_R$ of
$D_R$ is perfectoid, and the proof is complete.
\end{proof}

\begin{remark}\label{rem_tilt-of-site}
(i)\ \
Let $\underline X$ be any perfectoid quasi-affinoid scheme;
taking into account proposition \ref{prop_Scholze-tilt}(iii)
we deduce that the continuous map $\Spa(\bar u_{\underline X})$
of \eqref{subsec_what-name} induces a natural isomorphism of
sites (notation of \eqref{subsec_smaller-rational-site}) :
$$
\bE:\cR(\underline X)\isom\cR(\bE(\underline X))
\qquad
R\mapsto\Spa(\bar u_{\underline X})(R).
$$
(ii)\ \
Moreover, proposition \ref{prop_Letta-non-caduto}(iii)
says that for every $R\in\Ob(\cR(\underline X))$ the
sub-presheaf $h''_R$ of $h''_{\underline X}$ is represented
by a perfectoid quasi-affinoid scheme $\underline Y$
(notation of remark \ref{rem_yoneda-rationals}(i)),
and taking into account \eqref{eq_naturality-bE} together
with the equivalence of categories
$$
\underline X\tdu\mathsf{q.Afd.Sch}_\mathrm{perf}\isom
\bE(\underline X)\tdu\mathsf{q.Afd.Sch}_\mathrm{perf}
$$
described in \eqref{subsec_upgrade-E}, we deduce easily
that $\bE(\underline Y)$ represents the sub-presheaf
$h''_{\bE(R)}$ of $h''_{\bE(\underline X)}$ (details left to
the reader).
\end{remark}

\begin{theorem}\label{th_perfectoid-adic-space}
For every perfectoid quasi-affinoid scheme $\underline X$,
the following holds :
\begin{enumerate}
\item
The presheaves $\cO^\wedge_{\Spa\underline X}$,
$\cO^{\wedge\circ}_{\Spa\underline X}$ and
$\cO^{\wedge+}_{\Spa\underline X}$ are sheaves of topological
rings on the site $\cR(\underline X)$ (notation of
\eqref{subsec_smaller-rational-site}).
\item
For every rational subset $R$ of\/ $\Spa\,\underline X$ and
every finite covering $\fU$ of $R$ consisting of rational
subsets, the augmented alternating \v{C}ech complex
$C^\bullet_\mathrm{alt}(\fU,\cO^\wedge_{\Spa\,\underline X})$
has strict differentials, and its cohomology has the
discrete topology.
\end{enumerate}
\end{theorem}
\begin{proof}(i): By virtue of corollary \ref{cor_corcor}(i)
and proposition \ref{prop_general-case}(i), if
$\cO^\wedge_{\Spa\,\underline X}$ is a sheaf, the same holds for
$\cO^{\wedge+}_{\Spa\,\underline X}$ and $\cO^{\wedge\circ}_{\Spa\,\underline X}$,
so it suffices to check the assertion for $\cO^\wedge_{\Spa\,\underline X}$.
Now, say that $\underline X=(X,\cT_{\!X},A^+_X)$. Notice first that,
for every rational subset $R$ of $\Spa\,\underline X$, and every
morphism $\phi_{\underline Y/\underline X}:\underline Y\to\underline X$
representing the sub-presheaf $h_R$ of $h_{\underline X}$, the
quasi-affinoid scheme $\underline Y^\wedge$ is also perfectoid
(proposition \ref{prop_Letta-non-caduto}(iii)). Moreover,
lemmata \ref{lem_invariance-by-loc-hens}(iv) and
\ref{lem_rat-in-rat} and proposition
\ref{prop_Letta-non-caduto}(i) imply that
$\phi_{\underline Y/\underline X}$ induces equivalences of categories
\set\begin{equation}\label{eq_refer-in-later-th}
\cR(\underline Y^\wedge)\isom\cR(\underline Y)\isom\cR(\underline X)/R.
\end{equation}
Hence, let $A_X:=\cO_{\!X}(X)$ and $A:=A^\circ_X$, so that
$\sGamma^\circ(\underline X)=(A,A^+_X,X)$
(notation of \eqref{subsec_Gamma-circ}); according
to lemma \ref{lem_standard-coverings}, every open
covering of $\Spa\,\underline X$ can be refined by
the standard covering $R_\bullet:=(R_0,\dots,R_n)$
attached to a given sequence $a_\bullet:=(a_0,\dots,a_n)$
of elements of $A$ that generates an (open) ideal
$J\subset A$ such that
$\Spec\,A/J\subset\Spec\,A\setminus X$; for every
subset $\Lambda\subset\{0,\dots,n\}$ we define the
rational subset $R_\Lambda\subset R$ and the quasi-affinoid
ring
$\underline A_{X,\Lambda}:=(A_{X,\Lambda},A^+_{X,\Lambda},X_\Lambda)$
as in remark \ref{rem_extract-proj}(i), so that
$R_\Lambda=\Spa\,\underline A_{X,\Lambda}$ for every such
$\Lambda$. As in the proof of theorem
\ref{th_henselians-are-sheaves}, we are then reduced
to checking that the natural map
\set\begin{equation}\label{eq_albatross}
\cO^\wedge_{\Spa\,\underline X}(\Spa\,\underline X)\to
\Equal\Bigl(\prod_{i=0}^n\xymatrix{\cO^\wedge_{\Spa\,\underline X}(R_i)
\ar@<-.5ex>[r] \ar@<.5ex>[r] &}
\prod_{i,j=0}^n\cO^\wedge_{\Spa\,\underline X}(R_{\{i,j\}})\Bigr)
\end{equation}
is an isomorphism of topological rings. To this aim, set
$\bE:=\bE(A)$; by proposition \ref{prop_Scholze-tilt},
we may assume that there exists a sequence
$e_\bullet:=(e_0,\dots,e_n)$ of elements of $\bE$ such that
$a_i=\bar u_A(e_i)$ for $i=0,\dots,n$. Then, for every subset
$\Lambda\subset\{0,\dots,n\}$, we consider the subring
$$
C_{R,\Lambda}:=A\Bigl[\bar u_A(e^{1/p^k}_j)/\bar u_A(e^{1/p^k}_i)~|~
(j,i,k)\in\{0,\dots,n\}\times\Lambda\times\N\Bigr]
\subset A_{X,\Lambda}.
$$

\begin{claim}\label{cl_can-proceed}
For every $\Lambda\subset\{0,\dots,n\}$ the following holds :
\begin{enumerate}
\item
$C_{R,\Lambda}$ is a ring of definition of $A_{X,\Lambda}$.
\item
The natural map $\omega^\wedge_\Lambda:A^\wedge_{X,\Lambda}\to
\cO^\wedge_{\Spa\,\underline X}(R_\Lambda)$ of remark
\ref{rem_extract-proj}(ii) is an isomorphism of topological
$A$-algebras.
\end{enumerate}
\end{claim}
\begin{pfclaim} Both assertions are clear when
$\Lambda=\{i\}$, for any $i=0,\dots,n$, due to proposition
\ref{prop_Letta-non-caduto}(i) and claim \ref{cl_add-roots}.
We reduce to this case as follows. Let $r$ be the cardinality
of $\Lambda$, and set
$$
\Sigma:=\{e_{i_1}\cdots e_{i_r}~|~
0\leq i_1\leq i_2\leq\cdots\leq i_r\leq n\}.
$$
It is easily seen that the family $(\bar u_A(e)~|~e\in\Sigma)$
generates an open ideal of $A$ whose radical equals the radical
of $J$. Set also $e_\Lambda:=\prod_{i\in\Lambda}e_i\in\Sigma$, and
notice that
$$
R_\Lambda=
R_A\Bigl(\frac{\bar u_A(e)}{\bar u_A(e_\Lambda)}~|~e\in\Sigma\Bigr)
\qquad
C_{R,\Lambda}=A\Bigl[\bar u_A(e^{1/p^k})/\bar u_A(e^{1/p^k}_\Lambda)~|~
e\in\Sigma,\ k\in\N\Bigr].
$$
After replacing $e_0$ by $e_\Lambda$, and the system $(e_0,\dots,e_n)$
by $\Sigma$, we are therefore reduced to the case where $R=R_0$,
to which proposition \ref{prop_Letta-non-caduto}(i) and claim
\ref{cl_add-roots} apply.
\end{pfclaim}

In view of claim \ref{cl_can-proceed}, we may now proceed
as in remark \ref{rem_extract-proj}(iv) : we set
$$
B:=A[T_0^{1/p^\infty},\dots,T_n^{1/p^\infty}]
$$
and consider the homomorphism of $\N[1/p]$-graded $A$-algebras
$$
u:B\to A_X[Y^{1/p^\infty}]
\qquad
T_i^{1/p^k}\mapsto\bar u_A(e_i^{1/p^k})\cdot Y^{1/p^k}
\qquad
\text{for $i=0,\dots,n$ and every $k\in\N$}.
$$
Set also $S:=\Spec\,A^\circ_X$; according to
\eqref{subsec_functor-Q-graded}, the map
$u$ induces a well defined morphism of $S$-schemes
$$
\Proj\,(u):X\to\Proj\,B.
$$
Recall also that the inclusion map $A[T_0,\dots,T_n]\to B$
induces an affine morphism of schemes $\psi:\Proj\,B\to\P^n_S$.
Let $\P^n_S=\Omega_0\cup\cdots\cup\Omega_n$ be the standard
affine covering by complements of hyperplanes, as in remark
\ref{rem_extract-proj}(iv), set $\Omega'_i:=\psi^{-1}\Omega_i$
and $U_i:=(\Proj\,u)^{-1}\Omega'_i$ for $i=0,\dots,n$; there
follows an affine open covering of $\Proj\,B$ and an open
covering of $X$ :
$$
\Proj\,B=\Omega'_0\cup\cdots\cup\Omega'_n
\qquad
X=U_0\cup\cdots\cup U_n.
$$
Set also $U_\Lambda:=\bigcap_{i\in\Lambda}U_i$ for every non-empty
$\Lambda\subset\{0,\dots,n\}$; by remark \ref{rem_extract-proj}(i)
we have
\set\begin{equation}\label{eq_Us-too}
\cO_{\!X}(U_\Lambda)=A_{X,\Lambda}
\qquad
\text{for every non-empty $\Lambda\subset\{0,\dots,n\}$}.
\end{equation}
Let $V$ be the schematic image of $\Proj\,u$ (see remark
\ref{rem_extract-proj}(v)); a simple inspection shows
that
$$
V=\Proj\,\cA(e_\bullet)
$$
where $\cA(e_\bullet)$ is the $\N[1/p]$-graded algebra
associated with the sequence $e_\bullet$ via example
\ref{ex_E-of-beta} (cp. the proof of claim
\ref{cl_think-of-llamas}). Especially, we see that
$$
V_\Lambda:=V\cap\bigcap_{i\in\Lambda}\Omega'_i=
\Spec\,C_{R,\Lambda}
\qquad
\text{for every non-empty subset $\Lambda\subset\{0,\dots,n\}$}.
$$
Fix a finitely generated ideal $I$ of adic definition for
$\bE$, and for every $t\in\R_+$, denote by $C^\bullet_t$ the
augmented alternating \v{C}ech complex of
$I^{\lceil t\rceil}\cO_V$, relative to the open covering
$(V_i~|~i=0,\dots,n)$. Notice that if $\Lambda\neq\emptyset$,
the scheme $V_\Lambda$ is affine, hence we have
$I^{\lceil t\rceil}C_{R,\Lambda}=H^0(V_\Lambda,I^{\lceil t\rceil}\cO_V)$.
Moreover, by corollary \ref{cor_angular-blow-up} (and theorem
\ref{th_Cech-resolve}(ii)), the complex $C^\bullet_t$ is acyclic,
and the natural map $I^{\lceil t\rceil}A\to C^{-1}_t$ is an
isomorphism, for every $t\in\R_+$. It follows that
$$
I^{\lceil s\rceil}C^i_t=C^i_{t+s}
\qquad
\text{for every $i\in\Z$ and every $s,t\in\R_+$}
$$
and we endow $C^i_t$ with the $I$-adic topology, so that
$(C^i_{s+t}~|~s\in\R_+)$ is a fundamental system of open
submodules of $C^i_t$, for every $i\in\Z$ and every $t\in\R_+$.

Furthermore, denote by $A^\bullet_X$ the augmented
alternating \v{C}ech complex of $\cO_{\!X}$ relative to the
open covering $(U_i~|~i=0,\dots,n)$, and for every $i\in\Z$
endow $A^i_X$ with the topology deduced from the topologies
of the rings $A_{X,\Lambda}$, via the identifications
\eqref{eq_Us-too}, as in \eqref{subsec_smaller-rational-site}.
There follows a natural monomorphism of complexes of topological
$A$-modules
$$
C^\bullet_t\to A^\bullet_X
\qquad
\text{for every $t\in\R_+$}
$$
and the image of $C^i_t$ is open in $A^i_X$, for every
$i\in\Z$. Denote by $(A^{\wedge\bullet}_X,d_X^\bullet)$ and
$(C^{\wedge\bullet}_t,d^\bullet_t)$ the complexes obtained
by termwise completion of $A_X^\bullet$ and $C^\bullet_t$;
for every $t\in\R_+$ the induced map
$$
C^{\wedge\bullet}_t\to A^{\wedge\bullet}_X
$$
is still a monomorphism. By the foregoing, we have a natural
isomorphism of complexes
$$
C^{\wedge\bullet}_t\isom\lim_{s\in\R_+}C^\bullet_t/C^\bullet_{s+t}
\qquad
\text{for every $t\in\R_+$}
$$
and moreover $C^\bullet_t/C^\bullet_{s+t}$ is acyclic for every
$s,t\in\R_+$. By virtue of \cite[Th.3.5.8]{We}, it
follows that $C^{\wedge\bullet}_t$ is acyclic for every
$t\in\R_+$. Next, notice that the induced morphism
$$
A^\bullet_X/C^\bullet_t\to A^{\wedge\bullet}_X/C^{\wedge\bullet}_t
$$
is an isomorphism of complexes, for every $t\in\R_+$; we
conclude that
$$
H^nA^{\wedge\bullet}_X=
H^n(A^{\wedge\bullet}_X/C^{\wedge\bullet}_t)=
H^n(A^\bullet_X/C^\bullet_t)=H^nA^\bullet_X
\qquad
\text{for every $n\in\Z$}.
$$
On the other hand, claim \ref{cl_can-proceed}(ii) and
remark \ref{rem_extract-proj}(iii) yield an isomorphism
of complexes of topological $A$-modules
\set\begin{equation}\label{eq_nofamily}
A^{\wedge\bullet}_X\isom
C^\bullet_\mathrm{alt}(R_\bullet,\cO^\wedge_{\Spa\,\underline X}).
\end{equation}
Since $H^0A^\bullet_X=0$, we conclude that \eqref{eq_albatross}
is a ring isomorphism. Lastly, let $E$ be the equalizer
appearing in \eqref{eq_albatross}, and endow it with the
topology induced by the inclusion into
$\prod_{i=0}^n\cO^\wedge_{\Spa\,\underline X}(R_i)$; we have to
check that the resulting map
$\cO^\wedge_{\Spa\,\underline X}(\Spa\,\underline X)\to E$ is
open. In view of the isomorphism \eqref{eq_nofamily},
the assertion will follow from the following :

\begin{claim}
The differentials of the complex $(A^{\wedge\bullet}_X,d^\bullet_X)$
are strict and its cohomology has the discrete topology.
\end{claim}
\begin{pfclaim} For every $i\in\Z$, the family
$(\Ker\,d^i_t~|~t\in\R_+)$ yields a fundamental system of
open submodules of $\Ker\,d^i_X$. But we have already noticed
that $C^{\wedge\bullet}_t$ is acyclic, so this system is the
image of $(C^{\wedge i}_t~|~t\in\R_+)$, which is a fundamental
system of open subgroups of $A^{\wedge i}_X$, whence the contention.
\end{pfclaim}

(ii): We may argue as in the proof of theorem
\ref{th_compl-an-noeth-are-sheaves}(ii) : we choose
a standard covering $\fU'$ that refines $\fU$ and we
apply lemma \ref{lem_contract-cont}(i) and corollary
\ref{cor_double-strict-continuous} to the double complex
$C_\alt^\bullet(\fU,\fU')$ (details left to the reader).
\end{proof}

\sset\subsubsection{}\label{subsec_support-for-perf-qaff}
Let $\underline X:=(X,\cT_X,A^+_X)$ be any perfectoid
quasi-affinoid scheme; in light of theorem
\ref{th_perfectoid-adic-space}, we may now argue as
in \eqref{subsec_support-for-qaff} to first extend
$\cO^\wedge_{\Spa\,\underline X}$ to a sheaf of topological
rings on the topological space $\Spa\,\underline X$,
and then show that the stalks of the latter sheaf are
local rings. There follows a well defined morphism of
locally ringed spaces
$$
\sigma_{\underline X}:
(\Spa\,\underline X,\cO^\wedge_{\Spa\,\underline X})\to(X,\cO_{\!X})
$$
whose induced map on global sections is the identity map
$\cO_{\!X}(X)\to\cO^\wedge_{\Spa\,\underline X}(\Spa\,\underline X)$,
and whose underlying continuous map $\Spa\,\underline X\to X$
is the restriction of the support map of remark
\ref{rem_Spv-of-ring}(iii). Also corollary
\ref{cor_support-for-qaff} extends to the perfectoid case :

\begin{corollary}
With the notation of \eqref{subsec_support-for-perf-qaff},
the morphism $\sigma_{\underline X}$ induces an isomorphism
$$
H^i(X,\cO_{\!X})\isom
H^i(\Spa\,\underline X,\cO^\wedge_{\Spa\,\underline X})
\qquad
\text{for every $i\in\N$}.
$$
\end{corollary}
\begin{proof} We proceed as in the proof of corollary
\ref{cor_support-for-qaff} : first we show the following

\begin{claim}\label{cl_perf-loud}
In the situation of the corollary, suppose moreover that
$\underline X$ is affinoid. Then
$H^i(\Spa\,\underline X,\cO^\wedge_{\Spa\,\underline X})=0$ for
every $i>0$.
\end{claim}
\begin{pfclaim} Arguing as in the proof of claim \ref{cl_loud},
we reduce to checking that for every standard covering $\fU$
of $\Spa\,X$ the \v{C}ech cohomology
$H^i:=H^i_\mathrm{alt}(\fU,\cO^\wedge_{\Spa\,\underline X})$ vanishes
for every $i>0$; but the proof of theorem
\ref{th_perfectoid-adic-space}(i) shows that $H^i$ is
naturally isomorphic to the \v{C}ech cohomology
$H^i_\mathrm{alt}(U_\bullet,X)$, relative to a certain affine
open covering $U_\bullet$ of $X$. Then the assertion follows
from theorem \ref{th_Cech-resolve}.
\end{pfclaim}

Now, set $A_X:=\cO_{\!X}(X)$, say that
$X=\Spec\,A_X\setminus\Spec\,A_X/J$ for some finitely generated
ideal $J\subset A_X$, and pick a finite system of generators
$f_\bullet:=(f_0,\dots,f_n)$ for $J$. Let $R_\bullet:=(R_0,\dots,R_n)$
be the standard covering of $\Spa\,\underline X$ associated with
$f_\bullet$, and set as well $U_i:=\Spec\,A_X[f_i^{-1}]$ for
every $i=0,\dots,n$, so that $U_\bullet:=(U_0,\dots,U_n)$ is
an affine open covering of $X$, and moreover
$R_i=\sigma_{\underline X}^{-1}(U_i)$ for every $i=0,\dots,n$.
There follows a commutative diagram
$$
\xymatrix{ H^i_\alt(U_\bullet,\cO_{\!X}) \ar[r] \ar[d] &
H^i_\alt(R_\bullet,\cO^\wedge_{\Spa\,\underline X}) \ar[d] \\
H^i(X,\cO_{\!X}) \ar[r] &
H^i(\Spa\,\underline X,\cO^\wedge_{\Spa\,\underline X})
}$$
where the vertical arrows are isomorphisms, by virtue of
claim \ref{cl_perf-loud} and corollary \ref{cor_Leray}(ii).
We are then reduced to checking that the top horizontal
arrow is an isomorphism; but as already pointed out, the
latter assertion was shown in the proof of theorem
\ref{th_perfectoid-adic-space}(i).
\end{proof}

\sset\subsubsection{}\label{subsec_almost-vanish-O-plus}
Let again $\underline X$ be as in
\eqref{subsec_support-for-perf-qaff}; set $A_X:=\cO_{\!X}(X)$,
and let $V\subset X^+:=\Spec\,A^+_X$ be any quasi-compact
open subset containing the analytic locus of $X^+$, and
such that $\underline X$ spreads over $V$ (see definition
\ref{def_spread}(i)). Let also $J_V\subset A^+_X$ be the
unique radical ideal such that
$\Spec\,A^+_X/J_V=X^+\setminus V$. We consider the sheaf
of ideals
$$
\cJ_{\!V}\subset\cO^{\wedge+}_{\Spa\,\underline X}
\qquad
\text{on $\cR(\underline X)$}
$$
defined as the sheaf associated to the presheaf given
by the rule : $R\mapsto J_V\cO^{\wedge+}_{\Spa\,\underline X}(R)$
for every rational subset $R\subset\Spa\,\underline X$.

\begin{theorem}\label{th_integral-variant}
In the situation of \eqref{subsec_almost-vanish-O-plus}, we have :
\begin{enumerate}
\item
$\cJ_{\!V}(R)=J_V\cO^{\wedge+}_{\Spa\,\underline X}(R)$ for every
$R\in\Ob(\cR_V(\underline X))$ (notation of
\eqref{subsec_spread-site}).
\item
$H^i(X,\cJ_{\!V})=0$ for every $i>0$.
\end{enumerate}
\end{theorem}
\begin{proof} To begin with, we remark :

\begin{claim}\label{cl_J-transits-well}
Let $R\in\Ob(\cR_V(\underline X))$, and $\underline Y$
any topologically local quasi-affinoid scheme representing
the sub-presheaf $h_R$ of $h_{\Spa\,\underline X}$. Say that
$\underline Y^\wedge:=(Y^\wedge,\cT^\wedge_Y,A^{\wedge+}_Y)$.
We have :
\begin{enumerate}
\item
The completion $\underline Y^\wedge$ spreads over the
preimage $W$ of $V$ in $Y^{\wedge+}:=\Spec\,A^{\wedge+}_Y$.
\item
Let $J_W\subset A^{\wedge+}_Y$
be the unique radical ideal such that
$\Spec\,A^{\wedge+}_Y/J_W=Y^{\wedge+}\setminus W$. Then
$J_W=J_VA^{\wedge+}_Y$.
\end{enumerate}
\end{claim}
\begin{pfclaim}(i): Since $\underline Y$ spreads over $V$,
the assertion follows from proposition
\ref{prop_special-loci}(i).

(ii): We may find a finitely generated open ideal
$J_\bE\subset\bE^+$ such that $J_\bE^{\lceil 0\rceil}A^+_X=J_V$
(see \eqref{subsec_renzi-victor}); then we have
$\Spec\,A^{\wedge+}_Y/J_VA^{\wedge+}_Y=Y^{\wedge+}\setminus W$, and
on the other hand, arguing as in \eqref{subsec_renzi-victor}
we see that $J_\bE^{\lceil 0\rceil}A^{\wedge+}_Y$ is a radical ideal
of $A^{\wedge+}_Y$, whence the claim.
\end{pfclaim}

\begin{claim}\label{cl_terry-gilliam}
In order to prove the theorem it suffices to show that
for every $\underline X$ and $V$ as in
\eqref{subsec_almost-vanish-O-plus} we have
$\check{H}{}^0(\cR(\underline X),
J_V\cO^{\wedge+}_{\Spa\,\underline X})=J_V$ and
$\check{H}{}^i(\cR(\underline X),
J_V\cO^{\wedge+}_{\Spa\,\underline X})=0$ for every $i>0$.
\end{claim}
\begin{pfclaim} Indeed, in light of proposition
\ref{prop_change-site}(i,ii), both (i) and (ii) will follow
once we show that the presheaf $J_V\cO^{\wedge+}_{\Spa\,\underline X}$
is a sheaf on the site $\cR_V(\underline X)$, and
$H^i(\cR_V(\underline X),J_V\cO^{\wedge+}_{\Spa\,\underline X})=0$
for every $i>0$. Next, in view of theorem \ref{th_Cartan}, in
order to prove the latter assertions it suffices to show that
for every $R\in\Ob(\cR_V(\underline X))$ we have
$\check{H}{}^0(\cR_V(\underline X)/R,
J_V\cO^{\wedge+}_{\Spa\,\underline X})=J_V\cO^{\wedge+}_{\Spa\,\underline X}(R)$
and $\check{H}{}^i(\cR_V(\underline X)/R,
J_V\cO^{\wedge+}_{\Spa\,\underline X})=0$ for every
$i>0$. Proposition \ref{prop_change-site}(iii) then further
reduces to showing both identities, with $\cR_V(\underline X)$
replaced by $\cR(\underline X)$.
Lastly, for such a rational subset $R$, define $\underline Y$
and $W$ as in claim \ref{cl_J-transits-well}(i); in light
of claim \ref{cl_J-transits-well}(ii), and taking into
account the equivalences \eqref{eq_refer-in-later-th}, we may
then replace $\underline X$ by $\underline Y^\wedge$, and
$V$ by $W$, and reduce to prove the foregoing identities
for $R=\Spa\,\underline X$, whence the claim.
\end{pfclaim}

We consider first the case where the induced morphism
$\beta_{\underline X}:X\to X^+:=\Spec\,A^+_X$ is an open immersion,
so that $\Omega_{\underline X}$ is the image of $\beta_{\underline X}$
(notation of definition \ref{def_spread}), and we take
$V:=\Omega_{\underline X}$. By theorem
\ref{th_int-subrings-perfectoid}(iii) we have the perfectoid
quasi-affinoid ring
$$
\sGamma^+(\underline X):=(A^+_X,A^+_X,V)
$$
and clearly $\sSpec\,(\sGamma^+(\underline X))$ is isomorphic to
$\underline X$. In view of lemma \ref{lem_standard-coverings},
we deduce that every open covering of $\Spa\,\underline X$ can be
refined by the standard covering $\fU_\bullet:=(\fU_0,\dots,\fU_n)$
associated to some sequence $(a_0,\dots,a_n)$ of elements of
$A^+_X$ that generates an ideal $K$ with
$\Spec\,A^+_X/K\subset X^+\setminus V$, and notice that
$\fU_\bullet$ is a covering family of the site
$\cR_V(\underline X)$, due to proposition
\ref{prop_special-loci}(ii). We shall show more precisely
that the augmented alternating \v{C}ech complex
$C^\bullet_\mathrm{alt}(\fU_\bullet,J_V\cO^{\wedge+}_{\Spa\,\underline X})$
is acyclic. To this aim, set $\bE^+:=\bE(A^+_X)$; in view
of proposition \ref{prop_Scholze-tilt}, we may assume that
there exists a sequence $e_\bullet:=(e_0,\dots,e_n)$ of elements
of $\bE^+$ such that $a_i=\bar u_{A^+_X}(e_i)$ for $i=0,\dots,n$,
and we let $K_\bE\subset\bE^+$ be the open ideal generated by
$e_\bullet$ (proposition \ref{prop_Scholze-tilt}(ii)). Notice that
$J^{\lceil 0\rceil}_\bE$ and $K^{\lceil 0\rceil}_\bE$ are the radicals
of $J_\bE$ and respectively $K_\bE$; we deduce that
\set\begin{equation}\label{eq_Miereneuker}
J^{\lceil 0\rceil}_\bE=J^{\lceil 0\rceil}_\bE J^{\lceil 0\rceil}_\bE\subset
K^{\lceil 0\rceil}_\bE J^{\lceil 0\rceil}_\bE\subset J^{\lceil 0\rceil}_\bE.
\end{equation}
We consider the $\N[1/p]$-graded $A^+_X$-algebra $\cA$
associated as in example \ref{ex_E-of-beta} to the sequence
$e_\bullet$, {\em i.e.} such that
$$
\cA_0:=A^+_X
\qquad\text{and}\qquad
\cA_\gamma:=[e_\bullet]^{\La\gamma\Ra}A^+_X
\qquad
\text{for every $\gamma\in\N[1/p]\setminus\{0\}$}
$$
(with multiplication induced by that of $A^+_X$).
We have then a ring homomorphism as in
\eqref{subsec_angular-blow-up}
$$
\phi:\Z[T_0^{1/p^\infty},\dots,T^{1/p^\infty}_n]\to\cA
\qquad
T^{1/p^r}_i\mapsto\bar u_{A^+_X}(e_i^{1/p^r})
\qquad
\text{for $i=0,\dots,n$ and every $r\in\N$}
$$
whose image lies in $\cA_0=A^+_X$, whence corresponding
angular Rees algebras $\sR\La A^+_X\Ra$ and $\sR\La\cA\Ra$.
Set $Y_0:=\Proj\,\sR\La A^+_X\Ra$; the map $\phi$ induces a
morphism of schemes
$$
\Proj\,\phi:Y_0\to\P^n:=
\Proj\,\Z[T_0^{1/p^\infty},\dots,T^{1/p^\infty}_n].
$$
Let $\P^n=\Omega_0\cup\cdots\cup\Omega_n$ be the standard
affine open covering by complements of hyperplanes,
set $U_i:=(\Proj\,\phi)^{-1}\Omega_i$ for $i=0,\dots,n$,
and let
$$
C^\bullet:=C_\mathrm{alt}^\bullet(U_\bullet,\cO_{Y_0})
$$
be the  augmented alternating \v{C}ech complex relative
to the covering $U_\bullet:=(U_0,\dots,U_n)$. Taking into
account \eqref{eq_Miereneuker} and corollary
\ref{cor_ceiling-cohomol} (as well as theorem
\ref{th_Cech-resolve}(ii)), we deduce that
$$
C_\mathrm{alt}^\bullet(U_\bullet,J_V\cO_{Y_0})=J_VC^\bullet
$$
and moreover this second augmented alternating \v{C}ech
complex is acyclic. Furthermore, let $I_\bE$ be any
finitely generated ideal of adic definition for
$\bE^+$; corollary \ref{cor_angular-blow-up} shows
likewise that
$$
C_\mathrm{alt}^\bullet(U_\bullet,I^{\lceil t\rceil}_\bE\cO_{Y_0})
=I^{\lceil t\rceil}_\bE C^\bullet
\qquad
\text{for every $t\in\R_+$}
$$
and also these complexes are acyclic. We may describe
the terms of $C^\bullet$ as in the proof of theorem
\ref{th_perfectoid-adic-space} : for every
$i=0,\dots,n$, the $A^+_X$-module $C^i$ is a direct sum
$\bigoplus_{\Lambda\subset\{0,\dots,n\}}\cO_{Y_0}(U_\Lambda)$,
where $\Lambda$ ranges over the set of subsets of
cardinality equal to $i+1$, and
$U_\Lambda:=\bigcap_{j\in\Lambda}U_j$ for every such $\Lambda$.
Then
$$
\cO_{Y_0}(U_\Lambda)=
A^+_X[\bar u_{A^+_X}(e_i^{1/p^\infty})/\bar u_{A^+_X}(e_j^{1/p^\infty})~|~
(i,j)\in\{0,\dots,n\}\times\Lambda]\subset
A^+_X[a_j^{-1}~|~j\in\Lambda]
$$
for every non-empty $\Lambda\subset\{0,\dots,n\}$. If we
set likewise $\fU_\Lambda:=\bigcap_{j\in\Lambda}\fU_j$ for every
such $\Lambda$, we see that
$\cO_{Y_0}(U_\Lambda)\subset\cO^{\wedge+}_{\Spa\,\underline X}(\fU_\Lambda)$.
Set as usual $A_X:=\cO_{\!X}(X)$; in light of claim
\ref{cl_can-proceed} we see that there exists for
every $\Lambda\neq\emptyset$ a (unique) f-adic
topology on $A_X[a_j^{-1}~|~j\in\Lambda]$ for which the
localization map $A_X\to A_X[a_j^{-1}~|~j\in\Lambda]$
is f-adic and $\cO_{Y_0}(U_\Lambda)$ is a subring of definition.
Especially, the induced topology on $\cO_{Y_0}(U_\Lambda)$
is adic, and the natural map $A^+_X\to\cO_{Y_0}(U_\Lambda)$
is an adic ring homomorphism. Let us endow $C^i$ with the
corresponding product topology, for every $i=0,\dots,n$;
then $(I^{\lceil t\rceil}_\bE C^i~|~t\in\R_+)$ is a fundamental
system of open submodules of $C^i$, and hence also of $J_VC^i$,
for every such $i$. Moreover, we know that the natural maps
$J_VA^+_X\to J_VC^{-1}$ and
$I^{\lceil t\rceil}_\bE A^+_X\to I^{\lceil t\rceil}_\bE C^{-1}$ are
isomorphisms (again, by corollaries \ref{cor_ceiling-cohomol}(i)
and \ref{cor_angular-blow-up}(i)), so
$(I^{\lceil t\rceil}_\bE C^{-1}~|~t\in\R_+)$ is a fundamental
system of open submodules of $J_VC^{-1}$.
Taking into account \cite[Th.3.5.8]{We}, we conclude
that the termwise completion of $J_VC^\bullet$ is still acyclic.
Furthermore, the latter is isomorphic to $J_VC^{\wedge\bullet}$,
where $C^{\wedge\bullet}$ is the termwise completion of $C^\bullet$ :
indeed, $J_VC^{\wedge i}$ is dense and open in the completion of
$J_VC^i$ for every $i=0,\dots,n$, so the assertion is clear.

\begin{claim}\label{cl_not-yet-complete}
Set $W_\Lambda:=
U_\Lambda\cap\Spec\,\cO_{Y_0}(U_\Lambda)[a_j^{-1}~|~j\in\Lambda]$
for every non-empty $\Lambda\subset\{0,\dots,n\}$.
The datum $(\cO_{Y_0}(U_\Lambda),\cO_{Y_0}(U_\Lambda),W_\Lambda)$
is a quasi-affinoid ring, and its completion is a perfectoid
quasi-affinoid ring
$$
\underline R{}_\Lambda:=
(\cO_{Y_0}(U_\Lambda)^\wedge,\cO_{Y_0}(U_\Lambda)^\wedge,
W_\Lambda^\wedge:=
W_\Lambda\times_{U_\Lambda}\Spec\,\cO_{Y_0}(U_\Lambda)^\wedge).
$$
\end{claim}
\begin{pfclaim} Arguing as in the proof of claim
\ref{cl_can-proceed}, we may assume that
$\Lambda=\{i\}$, for some $i\in\{0,\dots,n\}$, in
which case it suffices to invoke proposition
\ref{prop_Letta-non-caduto}(iv).
\end{pfclaim}

With the notation of claim \ref{cl_not-yet-complete},
we see that
$$
\sGamma(\sSpec\,\underline R{}_\Lambda)=
(\cO_{Y_0}(U_\Lambda)^\wedge[a_j^{-1}~|~j\in\Lambda],
\cO_{Y_0}(U_\Lambda)^{\wedge+},W_\Lambda^\wedge)
$$
where $\cO_{Y_0}(U_\Lambda)^\wedge[a_j^{-1}~|~j\in\Lambda]$
carries the unique f-adic topology such that the
localization map $\cO_{Y_0}(U_\Lambda)^\wedge\to
\cO_{Y_0}(U_\Lambda)^\wedge[a_j^{-1}~|~j\in\Lambda]$ is open
(proposition \ref{prop_top-on-opens-fadic-case}(ii)),
and $\cO_{Y_0}(U_\Lambda)^{\wedge+}$ is the integral closure
of $\cO_{Y_0}(U_\Lambda)^\wedge$ in
$\cO_{Y_0}(U_\Lambda)^\wedge[a_j^{-1}~|~j\in\Lambda]$. In
other words, $\cO_{Y_0}(U_\Lambda)^\wedge[a_j^{-1}~|~j\in\Lambda]$
is the completion of
$\cO_{Y_0}(U_\Lambda)[a_j^{-1}~|~j\in\Lambda]=
A^+_X[a_j^{-1}~|~j\in\Lambda]$, where the latter is endowed
with the unique f-adic topology such that the localization
map $\cO_{Y_0}(U_\Lambda)\to A^+_X[a_j^{-1}~|~j\in\Lambda]$
is open (proposition \ref{prop_complete-f-adic}). It
follows that $\cO_{Y_0}(U_\Lambda)^{\wedge+}$ is also the
completion of the integral closure of $\cO_{Y_0}(U_\Lambda)$
in $A^+_X[a_j^{-1}~|~j\in\Lambda]$
(lemma \ref{lem_compl-and-int.clos}), which is the same
as the completion of the integral closure of
$A^+_X[a_i/a_j~|~(i,j)\in\{0,\dots,n\}\times\Lambda]$
in $A^+_X[a_j^{-1}~|~j\in\Lambda]$. This means that
$$
\sGamma(\sSpec\,\underline R{}_\Lambda)=
(\cO^\wedge_{\Spa\,\underline X}(\fU_\Lambda),
\cO^{\wedge+}_{\Spa\,\underline X}(\fU_\Lambda),W_\Lambda^\wedge).
$$
We have already observed that
$J_V\cO^{\wedge+}_{\Spa\,\underline X}(\fU_\Lambda)$
is a radical ideal of $\cO^{\wedge+}_{\Spa\,\underline X}(\fU_\Lambda)$;
combining with corollary \ref{cor_renzi-victor} we conclude that
$$
J_V\cO^{\wedge+}_{\Spa\,\underline X}(\fU_\Lambda)=J_V\cO_{Y_0}(U_\Lambda)^\wedge
\qquad
\text{for every non-empty $\Lambda\subset\{0,\dots,n\}$}
$$
which implies that the natural map of complexes
$J_VC^{\wedge\bullet}\to
C^\bullet_\mathrm{alt}(\fU_\bullet,J_V\cO^{\wedge+}_{\Spa\,\underline X})$
is an isomorphism, so the proof of the theorem is complete
in this case.

For the general case, let us set
$$
(X_\bE,\cT_\bE,\bE^+_X):=\bE(\underline X)
\qquad
X':=\beta^{-1}_{\underline X}(V)
\qquad
A_{X'}:=\cO_{\!X'}(X')
\qquad
\underline X':=X'\times_X\underline X
$$
(see example \ref{ex_restriction-of-qaff}).
Say that $\underline X'=(X',\cT_{X'},A^+_{X'})$;
then $\underline X'$ spreads over the preimage $V'$
of $V$ in $X'^+:=\Spec\,A^+_{X'}$ (proposition
\ref{prop_special-loci}(iv)); on the other hand, by
construction the image of $\beta_{\underline X'}:X'\to X'^+$
lies in $V'$, so $\beta_{\underline X'}$ is an open immersion
and $V'=\Omega_{\underline X'}$. Therefore we have as well
the perfectoid quasi-affinoid ring
$\sGamma(\underline X'):=(A^+_X,A^+_X,V')$, and
$\sSpec\,\sGamma^+(\underline X')$ is isomorphic to
$\underline X'$; especially, the latter is a perfectoid
quasi-affinoid scheme of the type considered in the
foregoing case. Next, since the natural morphism
$\pi:\underline X'\to\underline X$ is f-adic, proposition
\ref{prop_f-adic-Spectra}(iii.a) says that
$\Spa\,\pi:\Spa\,\underline X'\to\Spa\,\underline X$
induces a functor
$$
\cR(\underline X)\to\cR(\underline X')
\qquad
R\mapsto(\Spa\,\pi)^{-1}R
$$
which restricts to a functor
$$
\cR_V(\underline X)\to\cR_{V'}(\underline X').
$$
Indeed, let $R$ be any rational subset of
$\Spa\,\underline X$, and say that the corresponding
sub-presheaf $h_R$ of $h_{\Spa\,\underline X}$ is represented
by an f-adic morphism $\underline Y\to\underline X$ of
quasi-affinoid schemes. Then $(\Spa\,\pi)^{-1}R$ is
represented by the induced morphism
$\underline Y':=(\underline Y\times_{\underline X}\underline X')_\loc
\to\underline X'$ (see remark \ref{rem_yoneda-rationals}(iii)).
Now, if $\underline Y$ spreads over $V$, proposition
\ref{prop_special-loci}(i,iii,iv) implies that
$\underline Y'$ spreads over $V'$, whence the contention.
Furthermore, let $J_{V'}\subset A^+_{X'}$
be the unique radical ideal such that
$\Spec\,A^+_{X'}/J_{V'}=X'^+\setminus V'$; we have :

\begin{claim}\label{cl_looper}
For every $R\in\Ob(\cR_V(\underline X))$, the induced map
\set\begin{equation}\label{eq_looper}
J_V\cO^{\wedge+}_{\Spa\,\underline X}(R)\to
J_{V'}\cO^{\wedge+}_{\Spa\,\underline X'}((\Spa\,\pi)^{-1}R)
\end{equation}
is an isomorphism.
\end{claim}
\begin{pfclaim} We consider first the case where
$R=\Spa\,\underline X$ : arguing as in the proof of
claim \ref{cl_J-transits-well}(ii) we easily see that
$J_{V'}=J_VA^+_{X'}$. Therefore, let $\bar A{}^+_X$ be the
image of $A^+_X$ in $A^+_{X'}$; combining with corollary
\ref{cor_renzi-victor}, it follows that
$$
J_{V'}=J_V\bar A{}^+_X
$$
which already shows that the natural ring homomorphism
$\rho:A^+_X\to A^+_{X'}$ restricts to a surjection
$J_V\to J_{V'}$. Next, notice that the support of
$\Ker\,\rho$ lies in the closed subset $X^+\setminus V$,
and therefore $J_V\cap\Ker\,\rho$ is contained in the
nilradical of $A^+_X$; but $A^+_X$ is reduced (corollary
\ref{cor_perf-are-reduced}(i)), whence the claim, in this
case. Next, let $R\in\Ob(\cR_V(\underline X))$ be a general
rational subset; to ease notation, set
$R':=(\Spa\,\underline\pi)^{-1}R$, and let $\underline Z$
be any topologically local quasi-affinoid scheme representing
the sub-presheaf $h_{R'}$ of $h_{\Spa\,\underline X'}$. Define
$\underline Y^\wedge$, $Y^\wedge$, $Y^{\wedge+}$, $W\subset Y^{\wedge+}$
and $J_W$ as in claim \ref{cl_J-transits-well}. Let also
$Y':=\beta^{-1}_{\underline Y^\wedge}(W)$ and
$\underline Y':=Y'\times_{Y^\wedge}\underline Y^\wedge$.
Say that $\underline Y'=(Y',\cT_{Y'},A^+_{Y'})$, and let
moreover $W'$ be the preimage of $W$ in $Y'^+:=\Spec\,A^+_{Y'}$,
and $J_{W'}\subset A^+_{Y'}$ the unique radical ideal such that
$\Spec\,A^+_{Y'}/J_{W'}=Y'^+\setminus W'$; arguing as in the
foregoing, we see that $\underline Y'$ is a perfectoid
quasi-affinoid scheme that spreads over $W'$, and a
direct inspection of the constructions shows that
$\underline Y'$ is the completion of $\underline Z$.
Thus, the image of $\Spa\,\underline Y^\wedge$ (resp.
of $\Spa\,\underline Y'$) in $\Spa\,\underline X$
(resp. in $\Spa\,\underline X'$) agrees with $R$ (resp.
with $R'$), and taking into account claim
\ref{cl_J-transits-well}(ii) we conclude that
\eqref{eq_looper} is naturally identified with the map
$$
J_W\to J_{W'}
$$
induced by the projection
$\underline Y'\to\underline Y^\wedge$. By the previous case,
we know already that this map is bijective, whence the claim.
\end{pfclaim}

We may now conclude the proof of the theorem : indeed, we
have already remarked that every covering family of
$\cR(\underline X)$ can be refined by a standard covering
$\fU_\bullet:=(\fU_0,\dots,\fU_n)$ which is also a covering
family of $\cR_V(\underline X)$. It follows that
$\fU'_\bullet:=((\Spa\,\pi)^{-1}\fU_i~|~i=0,\dots,n)$ is a
standard covering of $\Spa\,\underline X'$ which is also
a covering family of $\cR_{V'}(\underline X')$. By the
previous case, we know already that
$\check{H}{}^0(\fU'_\bullet,J_{V'}\cO^{\wedge+}_{\Spa\,\underline X'})
=J_{V'}$ and $\check{H}{}^i(\fU'_\bullet,
J_{V'}\cO^{\wedge+}_{\Spa\,\underline X'})=0$ for every $i>0$.
In view of claim \ref{cl_looper}, it follows that
$\check{H}{}^0(\fU_\bullet,J_V\cO^{\wedge+}_{\Spa\,\underline X'})
=J_V$ and $\check{H}{}^i(\fU_\bullet,
J_V\cO^{\wedge+}_{\Spa\,\underline X'})=0$ for every $i>0$,
and taking into account claim \ref{cl_terry-gilliam},
the theorem follows.
\end{proof}

\subsection{Almost purity}\label{sec_almost-pure}
To begin with, let $p$ be a prime integer; we consider a
perfect $\F_p$-algebra $V$ and a radical ideal $\fm\subset V$.
Since $V$ is perfect, $(V,\fm)$ is clearly a basic setup,
in the sense of \cite[\S2.1.1]{Ga-Ra}. Our first observation
is the following :

\begin{lemma}\label{lem_lanzmann-dead}
In the situation of \eqref{sec_almost-pure}, let $A\subset B$
be an inclusion of perfect $V$-algebras, and $\fm_0\subset\fm$
a subideal whose radical is $\fm$. We set
$$
A':=\{x\in B~|~\fm\cdot x\subset A\}
\qquad
B':=\{x\in B~|~
\text{there exists $n\in\N$ with $\fm_0^n\cdot x\subset A$}\}.
$$
Then $A'$ and $B'$ are perfect $V$-algebras, and $A'$ is
integrally closed in $B'$.
\end{lemma}
\begin{proof} It is easily seen that $A'$ and $B'$ are
$V$-subalgebras of $B$, and obviously $A\subset A'\subset B'$.
Next, let $x\in B$, and say that $x^p\in A'$; thus,
$\fm\cdot x^p\subset A$, and therefore $(ax)^p\in A$ for
every $a\in\fm$. Since $A$ is perfect, it follows that
$ax\in A$ for every such $a$, so that $x\in A'$; this show
that $A'$ is perfect. Likewise one checks that $B'$ is perfect.

Next, let $x\in B'$, and suppose that $x^n+a_1x^{n-1}+\cdots+a_n=0$
for some $a_1,\dots,a_n\in A'$. Then
$A'[x]=A'+A'x+\cdots+A'x^{n-1}\subset B'$. By assumption,
there exists $k\in\N$ with $\fm_0^k\cdot x\subset A$, and
$\fm_0A'\subset\fm A'\subset A$; it follows easily that
$\fm_0^{nk}A'[x]\subset A$. Now, let $a\in\fm$; by assumption,
there exists $M\in\N$ such that $a^{p^M}\in\fm_0^{nk}$, whence
$(ax)^{p^M}\in A$, and since $A$ is perfect, we conclude that
$ax\in A$ for every such $a$, {\em i.e.} $x\in A'$. This
shows that $A'$ is integrally closed in $B'$.
\end{proof}

We may now state the following generalization of
\cite[Th.3.5.28]{Ga-Ra} :

\begin{theorem}\label{th_perfect-purity}
In the situation of \eqref{sec_almost-pure}, set $S:=\Spec\,V$
and $Z:=\Spec\,V/\fm$. Then :
\begin{enumerate}
\item
The pair $(S,Z)$ is almost pure, relative to the basic
setup $(V,\fm)$.
\item
If $Z$ is constructible in $S$, the pair $(S,Z)$ is normal,
for the basic setup $(V,\fm)$.
\end{enumerate}
\end{theorem}
\begin{proof}(ii): Denote by $j:U\to S$ the open immersion;
the kernel $\cK$ of the natural map $\cO_S\to j_*\cO_U$ is
the largest $\cO_S$-submodule such that $\fm\cK=0$, so the
induced morphism $\cO^a_S\to j_*\cO^a_U$ is a monomorphism.
It remains to check that the image of $\cO^a_S$ is integrally
closed in $j_*\cO^a_U$, and since $Z$ is constructible, we
are reduced to showing that $V^a$ is integrally closed in
$\cO_S(U)^a$ (\cite[Ch.I, Prop.9.2.1]{EGAI}). The latter follows
from lemma \ref{lem_lanzmann-dead} and \cite[Lemma 8.2.28]{Ga-Ra}.

(i): Let $(Z_\lambda~|~\lambda\in\Lambda)$ be the cofiltered
system of closed constructible subsets of $S$ containing $Z$,
and for every $\lambda\in\Lambda$ let $\fm_\lambda\subset\fm$
be the unique radical ideal such that
$\Spec\,V/\fm_\lambda=Z_\lambda$. Set also $U:=S\setminus Z$ and
$U_\lambda:=S\setminus Z_\lambda$ for every $\lambda\in\Lambda$.
Then $\fm_\lambda=\fm_\lambda^2$ and notice that both $\fm_\lambda$
(for every $\lambda\in\Lambda$) and $\fm$ trivially fulfill
condition $(\bB)$ of \cite[\S2.1.6]{Ga-Ra}, since $V$ is perfect.
We then get an essentially commutative diagram
$$
{\diagram (\cO_S,\fm\cO_S)^a\Et_\mathrm{afr}
\ar[r]^-\rho \ar[d] & \cO_U\Et \ar[d] \\
(\cO_S,\fm_\lambda\cO_S)^a\Et_\mathrm{afr}
\ar[r]^-{\rho_\lambda} & \cO_{U_\lambda}\Et
\enddiagram}
\qquad
\text{for every $\lambda\in\Lambda$}.
$$
However, on the one hand, proposition
\ref{prop_limit-of-alm-struct} and corollary
\ref{cor_limit-of-alm-struct}(ii,iv,v) imply that the
category $(\cO_S,\fm\cO_S)^a\Et_\mathrm{afr}$ is equivalent
to the $2$-limit of the system
$((\cO_S,\fm_\lambda\cO_S)^a\Et_\mathrm{afr}~|~\lambda\in\Lambda)$;
on the other hand, it is clear that $\cO_U\Et$
is equivalent to the $2$-limit of the system
$(\cO_{U_\lambda}\Et~|~\lambda\in\Lambda)$. Thus, we
are reduced to checking that $\rho_\lambda$ is an
equivalence for every $\lambda\in\Lambda$, and we may
therefore assume from start that $Z$ is constructible.
In view of (ii) and lemma \ref{lem_pure-almost-crit}(ii), we
are then reduced to showing that for every finite \'etale
$\cO_U$-algebra $\cB$ there exists an \'etale $\cO_S^a$-algebra
$\cA$ of almost finite rank such that $\rho(\cA)$ is isomorphic
to $\cB$. Now, by applying proposition \ref{prop_approx-al-alg}
with $X:=S$ and $\fm_0=\fm:=V$ (so, we consider here the
``classical limit'' almost ring structure) we obtain a
quasi-coherent $\cO_S$-algebra $\cR_0$ which is finitely
presented as an $\cO_S$-module, with an isomorphism
$\cB\isom\cR_{0|U}$ of $\cO_U$-algebras.

Let $\Phi_S:S\to S$ be the Frobenius endomorphism, and
consider the directed system of quasi-coherent
$\cO_S$-algebras $(\cR_n~|~n\in\N)$ with
$\cR_n:=\Phi_S^{-n}\cR_0$ for every $n\in\N$; explicitly,
$\cR_n(U')=\cR_0(U')$ for every open subset $U'\subset S$,
and the $\cO_S(U')$-algebra structure on $\cR_n(U')$ is
given by the composition of the Frobenius endomorphism
$\cO_S(U')\to\cO_S(U')$ with the structure morphism
$\cO_S(U')\to\cR_0(U')$. The transition maps
$\phi_n:\cR_n\to\cR_{n+1}$ are given on each such $U'$
by the Frobenius endomorphism $\cR_0(U')\to\cR_0(U')$,
for every $n\in\N$. We let
$$
\cR:=\colim_{n\in\N}\cR_n.
$$
Since $V$ is perfect, \cite[Th.3.5.13]{Ga-Ra} (applied
again in the ``classical limit'' case with $\fm=V$) implies
that $\phi_{n|U}$ is an isomorphism for every $n\in\N$,
so the natural morphism $\cR_0\to\cR$ induces an
isomorphism of $\cO_U$-algebras
$$
\cB\isom\cR_{|U}.
$$

\begin{claim}\label{cl_get-radical}
There exists a subideal $\fm_0\subset\fm$ whose radical
is $\fm$, and with
$$
\fm_0\cdot\Ker\,\phi_0=\fm_0\cdot\Coker\,\phi_0=0.
$$
\end{claim}
\begin{pfclaim} Since $Z$ is constructible, there exist
finitely many elements $x_1,\dots,x_n\in\fm$ such that
$U=\bigcup_{i=1}^n\Spec\,V[x_i^{-1}]$. Notice that both
$\cR_0$ and $\cR_1$ are finitely presented $\cO_S$-modules :
indeed, this holds for assumption for $\cR_0$, and then
it follows also for $\cR_1$, since $\Phi_S$ is an
automorphism of $S$. Thus, since moreover $\phi_{0|U}$
is an isomorphism, after replacing each $x_i$ by a power
$x_i^N$ for some suitable $N\in\N$, we may assume already
that $x_i\cdot\Coker\,\phi_0=0$ for every $i=1,\dots,n$.
Next, let $R_i:=\cR_i(S)$ for $i=0,1$, and denote by
$f:R_0\to R_1$ the Frobenius map; the natural map
$$
V[x_i^{-1}]\otimes_V\Hom_V(R_1,R_0)\to
\Hom_{V[x_i^{-1}]}(R_1[x_i^{-1}],R_0[x_i^{-1}])
$$
is bijective (\cite[Lemma 2.4.29(i.a)]{Ga-Ra}) and
$V[x_i^{-1}]\otimes_Vf$ is an isomorphism, hence we
may find for every $i=1,\dots,n$ a $V$-linear map
$g_i:R_1\to R_0$ such that
$V[x_i^{-1}]\otimes_V(g\circ f)=\one_{R_0}$. Then
there exists $N\in\N$ such that
$x_i^N\cdot(g\circ f)=x_i^N\cdot\one_{R_0}$ for every
$i=1,\dots,n$, whence $x_i^N\cdot\Ker\,f=0$ for
every $i=1,\dots,n$. The claim follows.
\end{pfclaim}

Let $\fm_0$ be as in claim \ref{cl_get-radical}, and for
every $n\in\N$ set $\fm_n:=\Phi_V^{-n}(\fm_0)$, where
$\Phi_V:V\to V$ denotes the Frobenius automorphism; then
$\fm_n\cdot\Ker\,\phi_n=\fm_n\cdot\Coker\,\phi_n=0$ for
every $n\in\N$. Notice that
$$
\fm^2_n\subset\fm_n\cdot\fm_{n+1}\cdots\fm_{n+k}
\qquad
\text{for every $n,k\in\N$}
$$
and therefore $\fm^2_n$ annihilates the kernel and cokernel
of the transition map $\cR_n\to\cR_{n+k}$, for every $n,k\in\N$;
we conclude that
\set\begin{equation}\label{eq_estimate}
\fm^2_n\cdot\Ker\,(\cR_n\to\cR)=
\fm^2_n\cdot\Coker\,(\cR_n\to\cR)=0
\qquad
\text{for every $n\in\N$}.
\end{equation}
Especially, $\cA:=\cR^a$ is an almost finitely presented
and uniformly almost finitely generated $\cO^a_S$-module
(\cite[Cor.2.3.13]{Ga-Ra}). Next, pick again elements
$x_1,\dots,x_n\in\fm$ such that $U$ is the union of the
affine open subsets $U_i:=\Spec\,V[x_i^{-1}]$ for $i=1,\dots,n$;
since $\cB$ is a projective $\cO_U$-module of finite rank,
for every $i=1,\dots,n$ there exist $N_i\in\N$ and
$\cO_U$-linear morphisms
$$
\cB_{|U_i}\xrightarrow{f_i}\cO^{\oplus N_i}_{U_i}
\xrightarrow{g_i}\cB_{|U_i}
\qquad\text{such that}\qquad
g_i\circ f_i=\one_{\cB|U_i}.
$$
Then, arguing as in the proof of claim \ref{cl_get-radical}
we find for every $i=1,\dots,n$, an integer $M_i\in\N$
and $\cO_S$-linear morphisms
$$
\cR_0\xrightarrow{f_{0,i}}\cO^{\oplus N_i}_S
\xrightarrow{g_{0,i}}\cR_0
\qquad\text{such that}\qquad
g_{0,i}\circ f_{0,i}=x_i^{M_i}\cdot\one_{\cR_0}
$$
whence, for every $k\in\N$, and every $i=1,\dots,n$, a pair
of $\cO_S$-linear morphisms
$$
\cR_k\xrightarrow{f_{k,i}}\cO^{\oplus N_i}_S
\xrightarrow{g_{k,i}}\cR_k
\qquad\text{such that}\qquad
g_{k,i}\circ f_{k,i}=x_i^{M_i/p^k}\cdot\one_{\cR_0}.
$$
It follows easily that $x_i^{M_i/p^k}\cdot\Tor^{\cO_S}_j(\cR_k,\cN)=0$
for every $j>0$, every $i=1,\dots,n$, every $k\in\N$ and every
quasi-coherent $\cO_S$-module $\cN$. After taking colimits,
we conclude that $\cA$ is a flat $\cO_S^a$-module; then $\cA$
is also an almost projective $\cO_S$-module of finite rank, in
light of \cite[Prop.2.4.18(ii) and Rem.4.3.10(i)]{Ga-Ra}.
Next, let $e$ be the diagonal idempotent of the \'etale
$\cO_U$-algebra $\cB$ (see remark \ref{rem_idemp-and-traces}(i)).
Pick a finite system of elements $\eps_1,\dots,\eps_k\in\fm$
such that $\fm$ is the radical of the ideal $\sum_{i=1}^kV\eps_i$;
we may then find $N\in\N$ and for every $i=1,\dots,k$ an element
$e_i\in R_0\otimes_VR_0$ such that
\set\begin{equation}\label{eq_approx-diag-idempotent}
(\eps_i^N\cdot e_i)_{|\Spec\,V[\eps^{-1}]}=(\eps^N_i\cdot e)_{|\Spec\,V[\eps^{-1}]}
\qquad\text{and}\qquad
\mu_0(e_i)=\eps^N_i
\end{equation}
where $\mu_n:\cR_n\otimes_{\cO_S}\cR_n\to\cR_n$ is the
multiplication map, for every $n\in\N$. Next, let $a_1,\dots,a_n$
be a finite system of generators for the $V$-module $R_0$; we may
find $N'\in\N$ such that
$$
\eps_i^{N'}\cdot(a_j\otimes 1-1\otimes a_j)\cdot e_i=0
\qquad
\text{for every $i=1,\dots,k$ and every $j=1,\dots,n$}
$$
and after replacing $e_i$ by $\eps_i^{N'}e_i$ and $N$ by $N+N'$ we may
therefore assume that \eqref{eq_approx-diag-idempotent} still holds,
and additionally :
$$
e_i\cdot\Ker\,\mu_0=0
\qquad
\text{for every $i=1,\dots,k$}.
$$
The same $e_1,\dots,e_k$ can also be regarded as elements
of $R_n\otimes_VR_n$, and we deduce that
\set\begin{equation}\label{eq_transfer-to-R_n}
\mu_n(e_i)=\eps_i^{N/p^n}
\qquad\text{and}\qquad
e_i\cdot\Ker\,\mu_n=0
\qquad
\text{for $i=1,\dots,k$ and every $n\in\N$}.
\end{equation}
Set $R:=\cR(S)$ and denote by $e_{n,i}\in R\otimes_VR$ the
image of $e_i$ under the natural morphism $\cR_n\to\cR$,
for $i=1,\dots,k$. Say that $\eps_i^{N''}\in\fm_0$ for
$i=1,\dots,k$ (where $\fm_0$ is as in claim
\ref{cl_get-radical}); then $\eps_i^{N''/p^n}\in\fm_n$ for
every $n\in\N$ and $i=1,\dots,k$. Recall that the kernel
of the multiplication map $\mu:\cR\otimes_{\cO_S}\cR\to\cR$
is generated, as a $\cR\otimes_{\cO_S}\cR$-module, by the
system $(a\otimes 1-1\otimes a~|~a\in R)$; in light of
\eqref{eq_estimate} and \eqref{eq_transfer-to-R_n} we
deduce that
$$
\eps_i^{2N''/p^n}\cdot e_{n,i}\cdot\Ker\,\mu=0
\qquad
\text{for $i=1,\dots,k$}
$$
where $\mu:\cR\otimes_{\cO_S}\cR\to\cR$ is the multiplication map.
After replacing $e_{n,i}$ by $\eps_i^{2N''/p^n}\cdot e_{n,i}$ and
$N$ by $N+2N''$, we then obtain $\cR\otimes_{\cO_S}\cR$-linear
morphisms
$$
\nu_{n,i}:\cR\to\cR\otimes_{\cO_S}\cR
\qquad
x\mapsto(x\otimes 1)\cdot e_{n,i}
\qquad
\text{for every $n\in\N$ and $i=1,\dots,k$}
$$
such that $\mu\circ\nu_{n,i}=\eps_i^{N/p^n}\cdot\one_\cR$.
Set $R:=\cR(S)$ for every $n\in\N$; it follows easily
that $\eps_i^{N/p^n}\cdot\Ext^j_{R\otimes_VR}(R,M)=0$ for every
$i=1,\dots,k$, every $n\in\N$, every $j>0$ and every
$R\otimes_VR$-module $M$, so $\cA$ is an almost projective
$\cA\otimes_{\cO_S}\cA$-module, and the proof is concluded.
\end{proof}

Notice next that, in the situation of \eqref{sec_almost-pure},
the ring $V$ is reduced, hence every flat almost finitely
presented $V^a$-module has almost finite rank (lemma
\ref{lem_reduced-or-integral}(ii)). We can generalize this
observation as follows :

\begin{proposition}\label{prop_always-afr}
Let $(V,\fm)$ be a basic setup as in \eqref{sec_almost-pure},
and consider ideals
$$
I\subset J\subset\fm
$$
with $J$ finitely generated. Let $P$ be any flat
almost finitely presented $(V/I)^a$-module. We have :
\begin{enumerate}
\item
$P$ has almost finite rank.
\item
If\/ $\Spec\,V/\fm$ is a constructible subset of\/
$S:=\Spec\,V$, then $P$ has finite rank.
\end{enumerate}
\end{proposition}
\begin{proof} In light of corollary
\ref{cor_limit-of-alm-struct}(iv) we may assume that
$\fm$ is the radical of a finitely generated subideal
$\fm_0$ containing $J$, in which case $\Spec\,V/\fm$ is
constructible in $S$, and it suffices to show assertion
(ii). Then we may even replace  $J$ by $\fm_0$, and
assume that $J^{\lceil 0\rceil}V=\fm$ (notation of remark
\ref{rem_Witt-are-f-adic}(i)).

\begin{claim} We may assume that $J$ lies in the Jacobson
radical of $V$.
\end{claim}
\begin{pfclaim} Indeed, let $x_1,\dots,x_k$ be a finite
system of generators of $J$; the induced map
$$
V\to(1+J)^{-1}V\times\prod_{i=1}^kV[x_i^{-1}]
$$
is faithfully flat, and clearly $V[x_i^{-1}]\otimes_VP$ is
a projective $V[x_i^{-1}]\otimes_VV/I$-module of finite rank
for every $i=1,\dots,k$, so we are reduced to checking that
$(1+J)^{-1}V\otimes_VP$ is an almost projective
$(1+J)^{-1}V\otimes_V(V/I)^a$-module of finite rank, whence
the claim.
\end{pfclaim}

Henceforth we assume that $J$ lies in the Jacobson radical
of $V$. Now, suppose that we have $\Lambda^r_{(V/J)^a}(P/JP)=0$
for some $r\in\N$; then $J\cdot\Lambda^r_{(V/I)^a}P=\Lambda^r_{(V/I)^a}P$,
and since $J$ is tight and $\Lambda^i_{(V/I)}P$ is almost
finitely generated, we deduce that $\Lambda^r_{(V/I)^a}P=0$,
by \cite[Lemma 5.1.7]{Ga-Ra}. Thus, we are reduced to
checking the assertion in case $I=J$. Next, notice that
there exists $n\in\N$ such that
$(J^{\La 1\Ra})^nV=J^{\La n\Ra}V\subset J$ (lemma
\ref{lem_mon-fract-powers}(ii.a,iv)), so $J^{\La 1\Ra}V$
is also a tight ideal contained in the Jacobson radical
of $V$, and by the same token we are reduced to checking
that $Q:=P/J^{\La 1\Ra}P$ is a $(V/J^{\La 1\Ra})^a$-module of
finite rank. To this aim, let us consider the $S$-scheme
$$
Y:=\Proj\,\sR
\qquad\text{where}\qquad
\sR:=\bigoplus_{\gamma\in\N[1/p]}J^{\La\gamma\Ra}V
$$
as in \eqref{subsec_angular-blow-up}, as well as the
closed immersion
$$
i:Y_0:=Y\times_S\Spec\,V/J^{\La 1\Ra}V=
\Proj\,\sR/J^{\La 1\Ra}\sR\to Y.
$$
Let also $\cT$ be the discrete topology on $V$, and notice
that $(V,\cT)$ is a perfectoid ring. Notice furthermore that
$(J^{\La 1\Ra}V)^a=(J^{\lceil 1\rceil}V)^a$; in light of example
\ref{subsec_angular-blow-up}(i), theorem
\ref{th_ceiling-cohomol}, and the short exact sequence
of quasi-coherent $\cO_Y$-modules
$$
0\to J^{\La 1\Ra}\cO_Y\to\cO_Y\to i_*\cO_{Y_0}\to 0
$$
we then deduce an isomorphism of $V^a$-algebras
$$
(V/J^{\La 1\Ra}V)^a\isom H^0(Y_0,\cO_{Y_0})^a
$$
(details left to the reader). Say that $J$ is generated
by a finite system $x_1,\dots,x_k$; then $Y$ admits a
closed immersion
$Y\to\P^{k-1}_S:=\Proj\,V[T_1^{1/p^\infty},\dots,T_k^{1/p^\infty}]$,
whence an affine covering $Y=U_1\dots,U_k$ consisting of
preimages of complements of hyperplanes of $\P^k_S$, as
usual. Set $U_{0,i}:=U_i\times_YY_0$ for $i=1,\dots,k$;
we get an injective ring homomorphism
$$
H_0(Y_0,\cO_{Y_0})\to A:=\prod_{i=1}^kH^0(U_{0,i},\cO_{Y_0})
$$
and on the other hand, natural identifications of
$V$-algebras
$$
A_i:=\gr_0\sR[x_i^{-1}]\isom H^0(U_i,\cO_Y)
\qquad
A_i/J^{\La 1\Ra}A_i\isom H^0(U_{0,i},\cO_{Y_0})
\qquad
i=1,\dots,k.
$$
Since $Q$ is a flat $(V/J^{\La 1\Ra}V)^a$-module, the same
holds for $\Lambda^r_{(V/J^{\La 1\Ra}V)^a}Q$, for every $r\in\N$;
hence the induced map
$$
\Lambda^r_{(V/J^{\La 1\Ra}V)^a}Q\to
A^a\otimes_{V^a}\Lambda^r_{(V/J^{\La 1\Ra}V)^a}Q
=\Lambda^r_{A^a}(A^a\otimes_{V^a}Q)
$$
is a monomorphism, and we are reduced to checking that
$A_i^a\otimes_{V^a}Q$ is an $(A_i/J^{\La 1\Ra}A_i)^a$-module
of finite rank for every $i=1,\dots,k$. Now, let
$(B_i,JB_i)$ be the henselization of the pair
$(A_i,JA_i)$. The ideal $\fm B_i$ is generated by the
system $(x_i^{1/p^n}~|~n\in\N)$ consisting of regular
elements of $B_i$; by \cite[Th.2.1.12(ii.b)]{Ga-Ra}
we deduce that $\fm B_i$ is a $B_i$-module of homological
dimension $\leq 1$. Lastly, since
$B_i/J^{\La 1\Ra}B_i=A_i/J^{\La 1\Ra}A_i$, we may apply
\cite[Th.5.5.7(i)]{Ga-Ra} to find a flat almost finitely
presented $B_i^a$-module $Q'$ with an isomorphism
$Q'/J^{\La 1\Ra}Q'\isom A^a_i\otimes_{V^a}Q$; since the
localization $B_i\to B_i[x_i^{-1}]$ is a monomorphism,
and since $B_i[x_i^{-1}]\otimes_{B_i}Q'$ is a projective
$B_i[x_i^{-1}]$-module of finite rank, the $B_i^a$-module
$Q'$ has finite rank as well. But then obviously the same
holds for the $(A_i/J^{\La 1\Ra}A)^a$-module $A^a_i\otimes_{V^a}Q$,
as required.
\end{proof}

\sset\subsubsection{}\label{subsec_compat-alm-structures}
Let $A$ be a perfectoid ring, $\bE:=\bE(A)$, and
$\underline\alpha\in W(\bE)$ a distinguished element in
the kernel of $u_A:W(\bE)\to A$. Recall that $\bar u_A:\bE\to A$
induces a ring isomorphism $\omega:\bE/\alpha_0\bE\isom A/pA$
(remark \ref{rem_nice-topology}(ii)); hence, for every ideal
$I\subset A$ containing $pA$, there exists a unique ideal
$I_\bE\subset\bE$ containing $\alpha_0\bE$ such that
$\omega(I_\bE/\alpha_0\bE)=I/pA$.

\begin{definition}\label{def_compat-alm-structures}
In the situation of \eqref{subsec_compat-alm-structures},
let $(A,\fm)$ and $(\bE,\fn)$ be two basic setups. We say
that $(A,\fm)$ and $(\bE,\fn)$ are {\em compatible} if there
exists an ideal of definition $I\subset A$ such that $\omega$
induces a bijection $(\fn+I_\bE)/I_\bE\isom(\fm+I)/I$.
\end{definition}

\begin{remark}\label{rem_compat-alm-structures}
With the notation of definition \ref{def_compat-alm-structures},
let $M,N$ be two $A/I$-modules, and $h:M\to N$ an $A$-linear
map. By assumption, $\bar u_A$ induces a ring isomorphism
$\omega_I:\bE/I_\bE\isom A/I$, and the compatibility of $(A,\fm)$
and $(\bE,\fn)$ means that $\omega_I$ is an isomorphism
$$
(\bE/I_\bE,(\fn+I_\bE)/I_\bE)\isom(A/I,(\fm+I)/I)
$$
in the category $\cB$ of basic setups as in \cite[\S3.5]{Ga-Ra}.
Hence $h^a$ is an isomorphism of $(A,\fm)^a$-modules if and
only if $\omega_I^*(h^a):\omega_I^*(M^a)\to\omega_I^*(N^a)$ is an
isomorphism of $(\bE,\fn)^a$-modules.
\end{remark}

\begin{lemma}\label{lem_compat-alm-structures}
In the situation of \eqref{subsec_compat-alm-structures},
consider two compatible basic setups $(A,\fm)$ and $(\bE,\fn)$,
and let also $f:B\to B'$ be any continuous map in $A\tdu\Perf$.
Then the morphism $f^a:B^a\to B'^a$ is an isomorphism in
$(A,\fm)^a\Alg$ if and only if $\bE(f)^a:\bE(B)^a\to\bE(B')^a$
is an isomorphism in $(\bE,\fn)^a\Alg$.
\end{lemma}
\begin{proof} By assumption, there exists an ideal of
definition $I\subset A$ such that $\omega$ induces a bijection
$(\fm+I)/I\isom(\fn+I_\bE)/I_\bE$, and arguing as in the proof
of lemma \ref{lem_was-third-cond}(i), we find a sequence
$\beta_1,\dots,\beta_k$ of elements of $\bE$ such that $I$ is
generated by $\bar u_A(\beta_1),\dots,\bar u_A(\beta_k)$; then,
after replacing each $\beta_i$ by $\beta_i^{1/p^n}$ for some large
$n\in\N$, we may also assume that $p\in I^t$, with $t:=k(p-1)+1$.
In light of lemma \ref{lem_before-name}(ii), it follows easily
that $I_\bE$ is generated by $\beta_1,\dots,\beta_k$, and
$\alpha_0\in I_\bE^t$. Denote by $\cT$ and $\cT'$ the $I$-adic
topologies of $B$ and $B'$; by lemma \ref{lem_fontaine} and
proposition \ref{prop_change-topol}, the topological rings
$(B,\cT)$ and $(B',\cT')$ are still perfectoid. Now, let
$\gr^\bullet B$ (resp. $\gr^\bullet B'$) be the graded ring
associated with the $I$-adic filtration on $B$ (resp. on $B'$);
likewise, let $\gr^\bullet\bE(B)$ (resp. $\gr^\bullet\bE(B')$) be the
graded ring associated with the $I_\bE$-adic filtration on $\bE(B)$
(resp. on $\bE(B')$). By claim \ref{cl_trickier-than-usual}(ii),
the maps $\bar u_B$ and $\bar u_{B'}$ induce graded ring
isomorphisms
$$
\theta:\gr^\bullet\bE(B)\isom\gr^\bullet B
\qquad
\theta':\gr^\bullet\bE(B')\isom\gr^\bullet B'.
$$
On the other hand, $f$ and $\bE(f)$ induce graded ring
homomorphisms
$$
\gr^\bullet(f):\gr^\bullet B\to\gr^\bullet B'
\qquad
\gr^\bullet\bE(f):\gr^\bullet\bE(B)\to\gr^\bullet\bE(B')
$$
such that :
\set\begin{equation}\label{eq_minute-detail}
\gr^\bullet(f)\circ\theta=\theta'\circ\gr^\bullet\bE(f).
\end{equation}
Now, if $f^a$ is an isomorphism, then the same holds for
$\gr^\bullet(f)^a$; conversely, if $\gr^\bullet(f)^a$ is an
isomorphism, then we easily see that
$f_n^a:=(A/I^n\otimes_Af)^a:(B/I^nB)^a\to(B'/I^nB')^a$ is an
isomorphism for every $n\in\N$, and therefore
$f^a=\lim_{n\in\N}f^a_n$ is an isomorphism (recall that the
functor $(-)^a$ commutes with limits, since it is a right
adjoint). Likewise we see that $\bE(f)^a$ is an isomorphism
if and only if the same holds for $\gr^\bullet\bE(f)^a$.
Thus, we are reduced to checking that $\gr^\bullet(f)^a$
is an isomorphism of $(A,\fm)^a$-algebras if and only if
$\gr^\bullet\bE(f)^a$ is an isomorphism of $(\bE,\fn)^a$-algebras.
However, notice that $\gr^\bullet B$ and $\gr^\bullet B'$ are
$A/I$-modules, and let $\omega_I:\bE/I_\bE\isom A/I$ be the
ring isomorphism deduced from $\bar u_A$; we may then
regard $\theta$ as an isomorphism
$\gr^\bullet\bE(B)\isom\omega_I^*(\gr^\bullet B)$ of
$\bE/I_\bE$-modules, and likewise for $\theta'$. Also,
$\gr^\bullet(f)$ is an isomorphism of $(A,\fm)^a$-algebras if
and only if $\omega^*_I(\gr^\bullet(f))$ is an isomorphism of
$(\bE,\fn)^a$-algebras, by remark \ref{rem_compat-alm-structures}.
Taking into account \eqref{eq_minute-detail}, the assertion follows.
\end{proof}

\sset\subsubsection{}\label{eq_down-to-business}
Let now $A$ be a perfectoid ring, and $\fm\subset A$ an open
radical ideal. It follows easily that $A^{\circ\circ}\subset\fm$;
moreover, if $\fm_\bE\subset\bE:=\bE(A)$ is associated with
$\fm$ as in \eqref{subsec_compat-alm-structures}, then
$\bar u_A$ induces a ring isomorphism $\bE/\fm_\bE\isom A/\fm$.

\begin{proposition}\label{prop_down-to-business}
In the situation of \eqref{eq_down-to-business}, let also
$f:B\to B'$ be any continuous map in $A\tdu\Perf$. The
following holds :
\begin{enumerate}
\item
$(A,\fm)$ and $(\bE,\fm_\bE)$ are basic setups fulfilling
condition {\em({\bf B})} of\/ \cite[\S2.1.6]{Ga-Ra}.
\item
More precisely, the multiplication maps
$\tilde\fm:=\fm\otimes_A\fm\to\fm$ and
$\tilde\fm_\bE:=\fm_\bE\otimes_\bE\fm_\bE\to\fm_\bE$
are isomorphisms.
\item
The morphism $f^a:B^a\to B'^a$ is an isomorphism in
$(A,\fm)^a\Alg$ if and only if\/ $\bE(f)^a:\bE(B)^a\to\bE(B')^a$
is an isomorphism in $(\bE,\fm_\bE)^a\Alg$.
\item
$\fm B$ is a radical ideal of $B$.
\end{enumerate}
\end{proposition}
\begin{proof}(i): Let $I$ be any ideal of definition of $A$.
Since $A/\fm$ is reduced, the same holds
for $\bE/\fm_\bE$, so $\fm_\bE$ is a radical ideal of $\bE$,
hence $\fm_\bE^2=\fm_\bE$, since $\bE$ is perfect. Now, let
$\underline\alpha\in\Ker\,u_A$ be any distinguished element;
by remark \ref{rem_nice-topology}(ii), the map $\bar u_A$
induces a ring isomorphism $\bE/\alpha_0\bE\isom A/pA$, and
it follows easily that $\fm^2+pA=\fm$. On the other hand,
we have $pA\subset I^2\subset\fm^2$, so $\fm^2=\fm$.
Next, it is clear that $\fm_\bE$ fulfills condition
({\bf B}), in particular $\fm_\bE$ is generated by its
$p$-th powers, therefore the same follows for the ideal
$\fm_\bE/\alpha_0\bE$ of $\bE/\alpha_0\bE$, and then
also for the ideal $\fm/pA$ of $A/pA$. Combining with
lemma \ref{lem_perfectoid}(iv), we deduce that $\fm$
is generated by its $p$-powers. If $l\neq p$ is any
prime integer, we have $l\in A^\times$, hence $\fm/l\fm=0$,
and taking into account \cite[Claim 2.1.9]{Ga-Ra}, it
follows that $\fm$ fulfills condition ({\bf B}).

(ii): Let us write $\fm$ as the union of a filtered
system $(\fm_\lambda~|~\lambda\in\Lambda)$ of open radical
subideals such that $\Spec\,A/\fm_\lambda$ is constructible
in $\Spec\,A$ for every $\lambda\in\Lambda$; clearly it
suffices to show that the corresponding map
$\fm_\lambda\otimes_A\fm_\lambda\to\fm_\lambda$ is an isomorphism
for every $\lambda\in\Lambda$. We may therefore assume that
$\fm$ is the radical of an open ideal generated by finitely
many elements $a_1,\dots,a_r$ of $A$; moreover, we may also
assume that there exist $\alpha_1,\dots,\alpha_r\in\bE$ with
$a_i=\bar u_A(\alpha_i)$ for $i=1,\dots,r$ (see remark
\ref{rem_nice-topology}(iii)). In this situation, we have
a ring homomorphism $u_0:R_{r,0}\to A$  such that
$u_0(T_i^{1/p^k})=\bar u_A(\alpha_i^{1/p^k})$ for $i=0,\dots,r$
(notation of \eqref{subsec_angular-blow-up}), and then it
is easily seen that $\fm=T^{\lceil 0\rceil}A$. Then the assertion
follows from proposition \ref{prop_back-to-toids} and theorem
\ref{th_ceiling-cohomol}(iv).

(iii): It is clear that the basic setups $(A,\fm)$ and
$(\bE,\fm_\bE)$ are compatible, in the sense of definition
\ref{def_compat-alm-structures}; then the assertion follows
from lemma \ref{lem_compat-alm-structures}.

(iv): The map $\bar u_B$ induces a ring isomorphism :
$$
\bE(B)/\fm_\bE\bE(B)\isom B/\fm B
$$
and we need to check that $B/\fm B$ is reduced; so we are
reduced to showing that the same holds
for $\bE(B)/\fm_\bE\bE(B)$, {\em i.e.} that $\fm_\bE\bE(B)$
is a radical ideal. But since $\bE(B)$ is perfect, this
is equivalent to proving that
$\Phi_{\bE(B)}(\fm_\bE B)=\fm_\bE B$ (where $\Phi_{\bE(B)}$
denotes the Frobenius automorphism of $\bE(B)$); the latter
is clear, since we have already shown that $\fm_\bE$ is a
radical ideal of $\bE$.
\end{proof}

\sset\subsubsection{}\label{subsec_almost-els-are-perfectoid}
In the situation of \eqref{eq_down-to-business}, let $B$
be any perfectoid $A$-algebra; according to proposition
\ref{prop_down-to-business}(i), we have a well defined
$(A,\fm)^a$-algebra $B^a$ and an $(\bE,\fm_\bE)^a$-algebra
$\bE(B)^a$. We denote by
$$
\eta_B:B\to B^a_*
\qquad\text{and}\qquad
\eta_{\bE(B)}:\bE(B)\to\bE(B)^a_*
$$
the units of adjunction, and endow $B^a_*$ (resp.
$\bE(B)^a_*$) with the unique ring topology such that
$\eta_B$ (resp. $\eta_{\bE(B)}$) is an adic morphism.

\begin{corollary}\label{cor_depressed}
With the notation of \eqref{subsec_almost-els-are-perfectoid},
the topological rings $B^a_*$ and $\bE(B)^a_*$ are perfectoid,
and there exists an isomorphism of\/ $\bE$-algebras
$$
\omega:\bE(B)^a_*\isom\bE(B^a_*)
\qquad\text{such that}\qquad
\omega\circ\eta_{\bE(B)}=\bE(\eta_B).
$$
\end{corollary}
\begin{proof} We consider first the case where the structure
map $A\to B$ is adic, and $Z:=\Spec\,A/\fm$ is constructible
in $X:=\Spec\,A$. Then, set $U:=\Spec\,B\times_X(X\setminus Z)$,
$B_U:=\cO_U(U)$, and let $\bar B\subset B_U$ be the image of the
restriction map $\rho:B\to B_U$. Notice that $B$ is reduced, by
corollary \ref{cor_perf-are-reduced}(i); it follows easily that
\set\begin{equation}\label{eq_rien-ne-vaut}
\fm\cdot\Ker\,\rho=0
\end{equation}
(the details are left to the reader). We remark :

\begin{claim}\label{cl_alm-els-of-B}
If $Z$ is constructible in $X$, we have a natural
identification of $B$-algebras :
$$
B^a_*\isom C:=\{x\in B_U~|~\fm\cdot x\subset\bar B\}.
$$
\end{claim}
\begin{pfclaim} Let $a_1,\dots,a_n\in A$ be a finite system of
elements that generate an ideal whose radical equals $\fm$; then :
$$
B_U=\xymatrix{
\Equal\,\bigl(\prod_{i=1}^nB[a_i^{-1}] \ar@<.5ex>[r] \ar@<-.5ex>[r] &
\prod_{i,j=1}^nB[a_i^{-1},a_j^{-1}]\bigr)}.
$$
Now, obviously $B[a_i^{-1}]^a_*=B[a_i^{-1}]$, and likewise
for $B[a_i^{-1},a_j^{-1}]$, for every $i,j=1,\dots,n$; since
the functor of almost elements is left exact, we deduce that
$B_{U*}^a=B_U$. Therefore, $\rho$ induces a map of $B$-algebras
$\rho^a_*:B^a_*\to B_U$, and notice that $\rho^a$ is a
monomorphism, due to \eqref{eq_rien-ne-vaut}, so the same holds
for $\rho^a_*$ (proposition \ref{prop_was-also-cofinal}(i)).
Obviously $\bar B{}^a=B^a$, hence
$\fm\cdot\rho^a_*(B^a_*)=\fm\bar B$, {\em i.e.} the image of
$\rho^a_*$ lies in the subring $C$. Moreover, $\rho$ induces
an isomorphism $C^a\isom B^a$, whence by adjunction, a unique
map of $B$-algebras $C\to B^a_*$ whose composition with
$\rho^a_*$ agrees with the inclusion map $C\to B_U$. This
shows that $\rho^a_*$ maps $B^a_*$ onto $C$, whence the claim.
\end{pfclaim}

Recall that $B_U$ carries a natural f-adic topology for
which $\bar B$ is a subring of definition (lemma
\ref{lem_topology-on-opens}(i)), and $\fm\bar B$ is
an open ideal of $\bar B$, since the structure map
$A\to B$ is adic; from claim \ref{cl_alm-els-of-B} it
then follows that $B^a_*$ is a bounded subring of $B_U$,
and therefore its topology as defined in
\eqref{subsec_almost-els-are-perfectoid} agrees with the
topology induced by the inclusion into $B_U$ (proposition
\ref{prop_f-adics}(i,ii)). From proposition
\ref{prop_general-case}(iii) we also see that
$(B^a_*)^{\circ\circ}=B_U^{\circ\circ}=B^{\circ\circ}$.

\begin{claim}\label{cl_perf-and-int-closed}
If the structure map $A\to B$ is adic, and
$Z$ is constructible in $X$, we have :

(i)\ \
$B^a_*/B^{\circ\circ}_U$ is a perfect $\F_p$-algebra and
it is integrally closed in $B_U^\circ/B^{\circ\circ}_U$.

(ii)\ \
$B^a_*$ is perfectoid and integrally closed in $B^\circ_U$.
\end{claim}
\begin{pfclaim}(i): Notice that $A/A^{\circ\circ}$,
$\bar B/B^{\circ\circ}$ and $B_U^\circ/B^{\circ\circ}_U$ are perfect
$\F_p$-algebras; also, $\fm/A^{\circ\circ}$ is a radical ideal of
$A/A^{\circ\circ}$, and since $Z$ is constructible, there exists
a finitely generated ideal $\fm_0\subset A/A^{\circ\circ}$ whose
radical is $\fm/A^{\circ\circ}$. With this notation,
$B_U^\circ/B^{\circ\circ}_U$ is an $A/A^{\circ\circ}$-algebra, and for
every $x\in B_U^\circ/B^{\circ\circ}_U$ there exists $k\in\N$ such
that $\fm_0^k\cdot x\subset\bar B/B^{\circ\circ}$. Then the
assertion follows from lemma \ref{lem_lanzmann-dead}.

(ii) follows from (i), claim \ref{cl_open-int-closed-subrings},
and theorem \ref{th_int-subrings-perfectoid}(iii).
\end{pfclaim}

Next, set $Z_\bE:=\Spec\,\bE/\fm_\bE$, $X_\bE:=\Spec\,\bE$, and
$U_\bE:=\Spec\,\bE(B)\times_{X_\bE}(X_\bE\setminus Z_\bE)$.
The foregoing argument applies also to $\bE(B)^a_*$,
and shows that the latter is a perfectoid subring of
$\bE_U:=\cO_{U_\bE}(U_\bE)$; moreover, we get cartesian
diagrams of topological rings :
$$
\xymatrix{ B^a_* \ar[r]^-{i_B} \ar[d] & B^\circ_U \ar[d] &
\bE(B)^a_* \ar[r]^-{\iota_\bE} \ar[d] &
\bE_U^\circ \ar[d] \\
B^a_*/B^{\circ\circ}_U \ar[r] & B^\circ_U/B^{\circ\circ}_U &
\bE(B)^a_*/\bE_U^{\circ\circ} \ar[r] & \bE^\circ_U/\bE_U^{\circ\circ}.
}$$
Furthermore, according to proposition \ref{prop_general-case}(iv),
the morphism of topological monoids
$\phi^{\flat\circ}_U:\bE^\circ_U\to B_U^\circ$ induces a ring
isomorphism
$\bar\phi^{\flat\circ}_U:\bE^\circ_U/\bE^{\circ\circ}_U\isom
B_U^\circ/B_U^{\circ\circ}$,
which in turn restricts to a ring isomorphism
$\bar\bE(B)/\bE_U^{\circ\circ}\isom\bar B/B^{\circ\circ}_U$,
where $\bar\bE(B)$ denotes the image of $\bE(B)$ into
$\bE_U$. Notice also that claim \ref{cl_alm-els-of-B}
yields natural identifications :
$$
\begin{aligned}
B^a_*/B^{\circ\circ}_U\isom&\,
\{x\in B_U^\circ/B_U^{\circ\circ}~|~
\fm\cdot x\subset\bar B/B^{\circ\circ}_U\} \\
\bE(B)^a_*/\bE_U^{\circ\circ}\isom&\,
\{x\in\bE_U^\circ/\bE_U^{\circ\circ}~|~
\fm_\bE\cdot x\subset\bar\bE(B)/\bE_U^{\circ\circ}\}.
\end{aligned}
$$
But clearly $\bar\phi^{\flat\circ}_U(\fm_\bE\cdot x)=
\fm\cdot\bar\phi^{\flat\circ}_U(x)$ for every
$x\in\bE_U^\circ/\bE_U^{\circ\circ}$, so $\bar\phi^{\flat\circ}_U$
restricts as well to a ring isomorphism
$$
\bE(B)^a_*/\bE_U^{\circ\circ}\isom B^a_*/B^{\circ\circ}_U.
$$
Summing up, and invoking proposition
\ref{prop_construct-new-perfs}, we obtain a unique
isomorphism of topological rings
$\omega:\bE(B)^a_*\isom\bE(B^a_*)$ fitting into the
commutative diagram :
$$
\xymatrix{
\bE(B)^a_* \ar[rr]^-{\iota_\bE}  \ar[d]_\omega & &
\bE^\circ_U \ar[d]^{\omega^\circ} \\
\bE(B^a_*) \ar[rr]^-{\bE(\iota_B)} & &
\bE(B^\circ_U)
}$$
where $\omega^\circ$ is the isomorphism provided
by theorem \ref{th_int-subrings-perfectoid}(i).
Lastly, denote by $\rho^\circ:B\to B^\circ_U$ and
$\rho^\circ_\bE:\bE(B)\to\bE_U^\circ$ the continuous
ring homomorphisms induced by the restriction maps;
the isomorphism $\omega^\circ$ is characterized by
the identity :
$$
\omega^\circ\circ\rho^\circ_\bE=\bE(\rho^\circ).
$$
Since $\rho^\circ$ and $\rho^\circ_\bE$ factor through
$\eta_B$ and respectively $\eta_{\bE(B)}$, we conclude
that $\omega\circ\eta_{\bE(B)}=\bE(\eta_B)$, as stated.
This completes the proof of the corollary in case
$Z$ is constructible.

Next, we consider the case where $Z$ can be an arbitrary
closed subset of $X$, and the structure map $A\to B$ is
still adic. Let us write $\fm$ as the filtered union of
a system $(\fm_\lambda~|~\lambda\in\Lambda)$ of open radical
subideals such that $\Spec\,A/\fm_\lambda$ is constructible
in $X$ for every $\lambda\in\Lambda$; correspondingly, we
obtain a filtered system
$(\fm_{\bE,\lambda}~|~\lambda\in\Lambda)$ of ideals in $\bE$
whose union is $\fm_\bE$. For every such $\lambda$ let
$B_\lambda$ (resp. $\bE(B)_\lambda$) be the image of $B$
(resp. of $\bE(B)$) in the category $(A,\fm_\lambda)^a\Alg$
(resp. $(\bE,\fm_{\bE,\lambda})^a\Alg$); we endow $B^a_{\lambda*}$
(resp. $\bE(B)^a_{\lambda*}$) with the unique adic topology
such that the unit of adjunction
$\eta_{B,\lambda}:B\to B^a_{\lambda*}$ (resp.
$\eta_{\bE(B),\lambda}:\bE(B)\to\bE(B)^a_{\lambda*}$) is adic.
After replacing $\Lambda$ by a cofinal subset, we may
assume that $\Lambda$ admits an initial element
$\lambda_0$. The localization functors
$(A,\fm)^a\Alg\to(A,\fm_\lambda)^a\Alg\to(A,\fm_\mu)^a\Alg$
induce natural maps of $B$-algebras
$$
f_\lambda:B^a_*\to B^a_{\lambda*}
\qquad\text{and}\qquad
f_{\lambda\mu}:B^a_{\lambda*}\to B^a_{\mu*}
\qquad
\text{for every $\lambda,\mu\in\Lambda$ with $\lambda\geq\mu$}
$$
which are clearly continuous. Moreover, clearly the images
$f^a_\lambda$ and $f^a_{\lambda\mu}$ of $f_\lambda$ and $f_{\lambda\mu}$
in $(A,\fm_{\lambda_0})^a\Alg$ are isomorphism, hence
$$
\fm_{\lambda_0}\cdot\Ker\,f_\lambda=
\fm_{\lambda_0}\cdot\Coker\,f_\lambda=0
\qquad
\fm_{\lambda_0}\cdot\Ker\,f_{\lambda\mu}=
\fm_{\lambda_0}\cdot\Coker\,f_{\lambda\mu}=0
\qquad
\text{for every $\lambda\geq\mu$}.
$$

\begin{claim}\label{cl_overlay-scrollbar}
The maps $f_\lambda$ and $f_{\lambda\mu}$ induce bijections
$$
\fm_{\lambda_0}B^a_*\isom\fm_{\lambda_0}B^a_{\lambda*}
\isom\fm_{\lambda_0}B^a_{\mu*}
\qquad
\text{for every $\lambda\geq\mu$}.
$$
\end{claim}
\begin{pfclaim} By the foregoing, we see already that
$f_\lambda$ (resp. $f_{\lambda\mu}$) maps $\fm_{\lambda_0}B^a_*$
onto $\fm_{\lambda_0}B^a_{\lambda*}$ (resp.
$\fm_{\lambda_0}B^a_{\lambda*}$ onto $\fm_{\lambda_0}B^a_{\mu*}$).
Moreover, by the previous case we also know that
$B^a_{\lambda*}$ is perfectoid, and under our current
assumptions $\fm_{\lambda_0}B^a_{\lambda*}$ is an open ideal
of $B^a_{\lambda*}$, hence $B^a_{\lambda*}/\fm_{\lambda_0}B^a_{\lambda*}=
B^a_{\lambda*}\hat\otimes_AA/\fm_{\lambda_0}$ is a discrete perfectoid
$\F_p$-algebra (proposition \ref{prop_stabil-cplted-tensors}(i));
especially, the latter is a reduced ring (corollary
\ref{cor_perf-are-reduced}(i)).
Now, let $x\in\fm_{\lambda_0}B^a_{\lambda*}\cap\Ker\,f_{\lambda\mu}$;
then $\fm_{\lambda_0}\cdot x=0$, whence $x^2=0$, and finally
$x=0$, which shows that $f_{\lambda\mu}$ is injective on
$\fm_{\lambda_0}B^a_{\lambda*}$ for every $\lambda\in\Lambda$.
Next, from corollary \ref{cor_limit-of-alm-struct}(i) we
know that $B^a_*$ is the limit of the cofiltered system
of $B$-algebras $(B^a_{\lambda*}~|~\lambda\in\Lambda)$,
hence the natural map
$$
\fm_{\lambda_0}B^a_*\to
L:=\lim_{\lambda\in\Lambda}\,\fm_{\lambda_0}B^a_{\lambda*}
$$
is injective. But we have also just seen that all the
projections $L\to\fm_{\lambda_0}B^a_{\lambda*}$ are isomorphisms,
so the maps $\fm_{\lambda_0}B^a_*\to\fm_{\lambda_0}B^a_{\lambda*}$
are injective as well.
\end{pfclaim}

From claim \ref{cl_overlay-scrollbar} and \cite[Lemma 3.5.3]{We}
we deduce that the natural map
$$
B^a_*/\fm_{\lambda_0}B^a_*\to
\lim_{\lambda\in\Lambda}\,B^a_{\lambda*}/\fm_{\lambda_0}B^a_{\lambda*}
$$
is a ring isomorphism; since we have already observed that
each quotient $B^a_{\lambda*}/\fm_{\lambda_0}B^a_{\lambda*}$ is
a perfect $\F_p$-algebra, it follows easily that the same
holds for $B^a_*/\fm_{\lambda_0}B^a_*$; especially, the latter
is a discrete perfectoid ring. We then get for every
$\lambda\in\Lambda$ a cartesian diagram of rings
$$
\cD_\lambda
\quad : \quad
{\diagram B^a_* \ar[r]^-{f_\lambda} \ar[d] &
B^a_{\lambda*} \ar[d] \\
B^a_*/\fm_{\lambda_0}B^a_* \ar[r] &
B^a_{\lambda*}/\fm_{\lambda_0}B^a_{\lambda*}
\enddiagram}
$$
which is also cartesian in the category of topological
rings, by lemma \ref{lem_adic-on-cartesian}, and proposition
\ref{prop_construct-new-perfs} implies that $B^a_*$ is
perfectoid. The same argument applies as well to the
perfectoid $\bE$-algebra $\bE(B)$, and it yields a
cartesian diagram of perfectoid rings :
$$
\cD_{\bE,\lambda}
\quad : \quad
{\diagram \bE(B)^a_* \ar[r]^-{f_{\bE,\lambda}} \ar[d] &
\bE(B)^a_{\lambda*} \ar[d] \\
\bE(B)^a_*/\fm_{\bE,\lambda_0}\bE(B)^a_* \ar[r] &
\bE(B)^a_{\lambda*}/\fm_{\bE,\lambda}\bE(B)^a_{\lambda*}.
\enddiagram}
$$
By the same token, the system $(f_{\bE,\lambda}~|~\lambda\in\Lambda)$
induces a ring isomorphism
$$
\bE(B)^a_*/\fm_{\bE,\lambda_0}\bE(B)^a_*\to\lim_{\lambda\in\Lambda}\,
\bE(B)^a_{\lambda*}/\fm_{\bE,\lambda_0}\bE(B)^a_{\lambda*}.
$$
Furthermore, by the previous case we already have an
isomorphism of topological rings
$$
\omega_\lambda:\bE(B)^a_{\lambda*}\isom\bE(B^a_{\lambda*})
\qquad\text{such that}\qquad
\omega_\lambda\circ\eta_{\bE(B),\lambda}=\bE(\eta_{B,\lambda})
\qquad
\text{for every $\lambda\in\Lambda$}.
$$
We may then identify naturally the bottom row of $\cD_\lambda$
with that of $\cD_{\bE,\lambda}$, and by invoking again proposition
\ref{prop_construct-new-perfs} we obtain a natural
isomorphism of topological rings
$\omega:\bE(B)^a_*\isom\bE(B^a_*)$ fitting into the
commutative diagrams
$$
{\diagram
\bE(B)^a_* \ar[rr]^-{f_{\bE,\lambda}} \ar[d]_\omega & &
\bE(B)^a_{\lambda*} \ar[d]^{\omega_\lambda} \\
\bE(B^a_*) \ar[rr]^-{\bE(f_\lambda)} & & \bE(B^a_{\lambda*})
\enddiagram}
\qquad
\text{for every $\lambda\in\Lambda$}.
$$
Then the sought identity $\omega\circ\eta_{\bE(B)}=\bE(\eta_B)$
follows from the corresponding ones for the isomorphisms
$\omega_\lambda$, after taking limits (details left to the
reader). This completes the proof, in case the topologies
of $A$ and $B$ are $p$-adic.

For the general case, let $\cT_{A,p}$ (resp. $\cT_{B,p}$) be
the $p$-adic topology on $A$ (resp. on $B$); then $(A,\cT_{A,p})$
and $(B,\cT_{B,p})$ are perfectoid topological rings as well
(proposition \ref{prop_change-topol}(i)), and the rings
underlying $\bE((B,\cT_{B,p})^a_*)$ and $\bE(B,\cT_{B,p})^a_*$
coincide with the rings underlying $\bE(B^a_*)$ and
respectively $\bE(B)^a_*$, so the previous case already yields
a ring isomorphism $\omega$ as sought. It remains to check
that $B^a_*$ (resp. $\bE(B)^a_*$) is perfectoid for the
topology such that $\eta_B$ (resp. $\eta_{\bE(B)}$) is adic,
and that $\omega$ is also an isomorphism of topological
rings for these topologies. However, let $I\subset B$ be
any finitely generated ideal of definition; the natural
map $B^a\to\lim_{n\in\N}(B/I^n)^a$ is an isomorphism, hence
the same holds for the induced map
$B^a_*\to\lim_{n\in\N}(B/I^n)^a_*$, and especially, $B^a_*$ is
complete and separated for the linear topology for which
$(J_n:=\Ker\,B^a_*\to(B/I^n)^a_*~|~n\in\N)$ is a fundamental
system of open submodules (corollary
\ref{cor_limits-and-complete}(i)). Clearly $I^n B^a_*\subset J_n$
for every $n\in\N$, so $B^a_*$ is also complete and separated
for its $I$-adic topology (lemma \ref{lem_fontaine}). Together
with proposition \ref{prop_change-topol}(ii), this proves that
$B^a_*$ is perfectoid, and the same applies to $\bE(B)^a_*$.
Lastly, since $\eta_B$ is adic, the same holds for $\bE(\eta_B)$
(corollary \ref{prop_change-topol}), but then the identity
$\omega\circ\eta_{\bE(B)}=\bE(\eta_B)$ implies easily that
$\omega$ is an isomorphism of topological rings.
\end{proof}

\sset\subsubsection{}\label{subsec_double-shriek-perf}
In the situation of \eqref{subsec_almost-els-are-perfectoid},
denote also by $\eps_B:B^a_{!!}\to B$ and
$\eps_{\bE(B)}:\bE(B)^a_{!!}\to\bE(B)$ the counits of adjunction,
and endow $B^a_{!!}$ (resp. $\bE(B)^a_{!!}$) with the unique
topology such that the natural map $A\to B^a_{!!}$ (resp.
$\bE\to\bE(B)^a_{!!}$) is adic.

\begin{corollary} With the notation of
\eqref{subsec_double-shriek-perf}, the topological rings
$B^a_{!!}$ and $\bE(B)^a_{!!}$ are perfectoid and there exists
an isomorphism of\/ $\bE$-algebras
$$
\omega:\bE(B)^a_{!!}\isom\bE(B^a_{!!})
\qquad
\text{such that}\qquad
\bE(\eps_B)\circ\omega=\eps_{\bE(B)}.
$$
\end{corollary}
\begin{proof} Notice that, by virtue of proposition
\ref{prop_down-to-business}(ii) we have cartesian
diagrams of rings
$$
\cD
\quad :\quad
{\diagram B^a_{!!} \ar[r]^-{\eps_B} \ar[d] & B \ar[d] \\
A/\fm \ar[r] & B/\fm B.
\enddiagram}
\qquad\qquad
\cD_\bE
\quad :\quad
{\diagram \bE(B)^a_{!!} \ar[r]^-{\eps_{\bE(B)}} \ar[d] &
\bE(B) \ar[d] \\
\bE/\fm_\bE \ar[r] & \bE(B)/\fm_\bE\bE(B).
\enddiagram}
$$
Moreover, clearly we may assume that the structure map
$A\to B$ is an adic ring homomorphism, in which case the
same holds for the map $\bE\to\bE(B)$ (theorem
\ref{th_adic-to-adic}(i)), and the bottom rows of $\cD$ and
$\cD_\bE$ consist of discrete perfectoid $\F_p$-algebras
(cp. the proof of claim \ref{cl_overlay-scrollbar}).
Then $\cD$ and $\cD_\bE$ are cartesian in the category of
topological rings as well (lemma \ref{lem_adic-on-cartesian}),
and the assertion follows from proposition
\ref{prop_construct-new-perfs}.
\end{proof}

\sset\subsubsection{}\label{subsec_formal-alm-flat}
 In the situation of \eqref{eq_down-to-business},
let $f:B\to C$ be a morphism of $A$-perfectoid algebras,
$(b_\lambda~|~\lambda\in\Lambda)$ a system of elements of
$B$, and for every $n\in\N$ denote by $I_n\subset B$ the
ideal generated by the system $(b^n_\lambda~|~\lambda\in\Lambda)$.

\begin{lemma}\label{lem_depressed}
With the notation of \eqref{subsec_formal-alm-flat},
suppose that $p^kB\subset I_1$ for some $k\in\N$, and
that $(f\otimes_BB/I_1)^a:(B/I_1)^a\to(C/I_1C)^a$ is a
flat morphism of $(A,\fm)^a$-algebras. Then the same
holds for $(f\otimes_BB/I_n)^a$, for every $n\in\N$.
\end{lemma}
\begin{proof} Set $J_n:=p^nB+I_n$ for every $n\in\N$;
since $p^kB\subset I_1$, it is easily seen that for
every $n\in\N$ there exists $n'\in\N$ such that
$J_{n'}\subset I_n\subset J_n$ (details left to the
reader). Hence, it suffices to show that
$(f\otimes_BB/J_n)^a$ is flat for every $n\in\N$. Next,
let $\underline\alpha\in\Ker\,u_B$ be a distinguished
element, and pick for every $\lambda\in\Lambda$
elements $e_\lambda\in\bE(B)$ and $c_\lambda\in B$ such that
$$
b_\lambda=\bar u_B(e_\lambda)+\bar u_B(\alpha_0)\cdot c_\lambda.
$$
For every $n\in\N$, denote by $I'_n\subset B$ the ideal
generated by $(\bar u_B(e_\lambda^n)~|~\lambda\in\Lambda)$,
and set $J'_n:=\bar u_B(\alpha_0^n)\cdot B+I'_n$; it follows
easily that
$$
J_{2n-1}\subset J'_n
\qquad\text{and}\qquad
J'_{2n-1}\subset J_n
\qquad
\text{for every $n\in\N$}
$$
so $(f\otimes_BB/J'_1)^a$ is flat, and we are further reduced
to checking that $(f\otimes_BB/J'_n)^a$ is flat for every
$n\in\N$. For every $n\in\N$, denote by  $I_{\bE,n}\subset\bE(B)$
the ideal generated by $(e^n_\lambda~|~\lambda\in\Lambda)$, and set
$J_{\bE,n}:=\alpha^n_0\bE(B)+I_{\bE,n}$; notice the natural
identifications
\set\begin{equation}\label{eq_depressed}
B/J'_1=\bE(B)/J_{\bE,1}
\qquad\text{and}\qquad
C/J'_1C=\bE(C)/J_{\bE,1}\bE(C).
\end{equation}
Denote also by $\Phi_{B/J'_1}$ and $\Phi_{C/J'_1C}$ the
Frobenius endomorphism of $B/J'_1$ and respectively
$C/J'_1C$. We remark :

\begin{claim}\label{cl_isabel}
The induced diagram of $\bE$-algebras
$$
\xymatrix{ B/J'_1 \ar[rr]^-{f\otimes_BB/J'_1}
\ar[d]_{\Phi_{B/J'_1}} & & C/J'_1C \ar[d]^{\Phi_{C/J'_1C}} \\
B/J'_1 \ar[rr]^-{f\otimes_BB/J'_1} & & C/J'_1C
}$$
is cocartesian.
\end{claim}
\begin{pfclaim} Indeed, since $\bE(B)$ is perfect,
it is easily seen that $\Phi_{B/J'_1}$ is surjective,
and its kernel is the ideal generated by the system
$\{\alpha^{1/p}_0\}\cup\{e^{1/p}_\lambda~|~\lambda\in\Lambda\}$,
under the identifications \eqref{eq_depressed};
likewise, $\Phi_{C/J'_1C}$ is surjective and its kernel
is generated by the image in $C/J'_1C$ of the same
system, whence the claim.
\end{pfclaim}

For every $n\in\N$, set
$\fn_n:=\{\underline a\in W_n\bE~|~a_0,\dots,a_{n-1}\in\fm\}$;
since $\fm_\bE$ fulfills condition $(\bB)$ of
\cite[\S2.1.6]{Ga-Ra}, proposition \ref{prop_alm-str-on-Witt}
says that $\fn_n$ is the unique ideal
$\fn_n\subset W_n\bE$ with $\fn^2_n=\fn_n$ and
$\bar\omega_i(\fn_n)=\fm_\bE$ for $i=0,\dots,n-1$.
From claim \ref{cl_isabel} and theorem
\ref{th_Witt-and-etale-maps}(i) it follows that the
induced morphism of $(W_n\bE,\fn_n)^a$-algebras
$$
W_n(f\otimes_BB/J'_1)^a:
W_n(\bE(B)/J_{\bE,1})^a\to W_n(\bE(C)/J_{\bE,1}\bE(C))^a
$$
is flat for every $n\in\N$. Moreover, the Frobenius
automorphisms of $\bE(B)$ and $\bE(C)$ induce isomorphisms
$\bE(B)/J_{\bE,1}\isom\bE(B)/J_{\bE,p^k}$ and
$\bE(C)/J_{\bE,1}\bE(C)\isom\bE(C)/J_{\bE,p^k}\bE(C)$ for every
$k\in\N$, whence a commutative diagram of $W_n\bE$-algebras
$$
{\diagram
W_n(\bE(B)/J_{\bE,1})
\ar[rrrr]^-{W_n(f\otimes_B\bE(B)/J_{\bE,1})}
\ar[d] & & & & W_n(\bE(C)/J_{\bE,1}\bE(C)) \ar[d] \\
W_n(\bE(B)/J_{\bE,p^k})
\ar[rrrr]^-{W_n(f\otimes_B\bE(B)/J_{\bE,p^k})} & & & &
W_n(\bE(C)/J_{\bE,p^k}\bE(C))
\enddiagram}
\quad
\text{for every $n,k\in\N$}
$$
whose vertical arrows are isomorphisms, and where the
rings on the bottom row are regarded as
$W_n\bE$-algebras via restriction of scalars along
the automorphism $W_n(\Phi_\bE):W_n\bE\to W_n\bE$.
Since we have
$$
W_n(\Phi_\bE)(\fn_n)=\fn_n
\qquad
\text{for every $n\in\N$}
$$
it follows that $W_n(f\otimes_B\bE(B)/J_{\bE,p^k})^a$
is a flat morphism of $(W_n\bE,\fn_n)^a$-algebras,
for every $n,k\in\N$. We remark, quite generally :

\begin{claim}\label{cl_instated}
Let $R$ be any ring such that the Frobenius endomorphism
$\Phi_{R/pR}$ of $R/pR$ is surjective,
$(x_\lambda~|~\lambda\in\Lambda)$ any set of elements of
$R$, and $I\subset R$ the ideal generated by this set.
Then, for every $n\in\N$, the kernel of the projection
$\pi_n:W_nR\to W_nR/I$ is generated by
$$
(V^i_R(x_\lambda)~|~\lambda\in\Lambda,\ i=0,\dots,n-1).
$$
\end{claim}
\begin{pfclaim} Arguing by induction on $n\in\N$, we
are reduced to checking that the ideal
$V_nI:=V_nR\cap\Ker\,\pi_{n+1}$ is generated by
$(V^nx_\lambda~|~\lambda\in\Lambda)$, for every $n\in\N$.
However, according to claim \ref{cl_W_n-structure} we
have an isomorphism of $W_{n+1}R$-modules : $I\isom V_nI$,
where $I$ is regarded as a $W_{n+1}R$-module via restriction
of scalars along the ring homomorphism
$\bar\omega_n:W_{n+1}R\to R$. Thus, we come down to showing
that $\bar\omega_n$ is a surjection for every $n\in\N$.
We argue by induction on $n$ : the assertion is trivial
for $n=0$, so suppose that $\bar\omega_n$ is already known
to be surjective for some $n\geq 0$; since $\Phi_{R/pR}$ is
surjective, the same holds for $\Phi^n_{R/pR}$, and taking
into account \eqref{eq_witt-pols}  together with our inductive
assumption, we easily deduce that $\bar\omega_{n+1}$ is surjective
as well.
\end{pfclaim}

In light of claim \ref{cl_instated} we see that
$$
W_n(\bE(B)/J_{\bE,p^k})=W(\bE(B))/(p^nW(\bE(B))+\cI_k)
$$
where $\cI_k$ is the ideal generated by
$(p^i\cdot\tau_{\bE(B)}(e^{p^{k-i}}_\lambda)~|~
\lambda\in\Lambda,\ i=0,\dots,k)$. Taking into account
lemma \ref{lem_drop-conditions}(iv), we deduce that
$$
A\otimes_{W(\bE)}W_n(\bE(B)/J_{\bE,p^k})=
B/(p^nB+I''_k)
\qquad
\text{for every $n,k\in\N$}
$$
where $I''_k$ is the ideal generated by
$(p^i\cdot\bar u_B(e_\lambda^{p^{k-i}})~|~
\lambda\in\Lambda,\ i=0,\dots,k)$. Notice that
\set\begin{equation}\label{eq_khazar-crap}
J'_{p^{2k-1}}\subset I''_{2k-1}\subset J'_k
\qquad
\text{for every $k\in\N$}.
\end{equation}
We have as well
$A\otimes_{W(\bE)}W_n(f\otimes_{\bE(B)}\bE(B)/J_{\bE,p^k})=
f\otimes_BB/(p^nB+I''_k)$, and furthermore, the ring homomorphism
$W_n\bE\to A/p^nA$ induced by $u_A$ maps $\fn_n$ onto
$\fm/p^nA$; therefore, $(f\otimes_BB/(p^nB+I''_k))^a$ is a
flat morphism of $(A,\fm)^a$-algebras, for every $n,k\in\N$.
Combining with \eqref{eq_khazar-crap}, the assertion follows.
\end{proof}

\begin{proposition}\label{prop_alm-etale-over-perf-is-perf}
In the situation of \eqref{eq_down-to-business}, let $B$ be
a perfectoid $A$-algebra, $f^a:B^a\to C$ an \'etale morphism
of $(A,\fm)^a$-algebras such that $C$ is almost finitely
presented as a $B^a$-module, and endow $C_*$ with the unique
ring topology such that the induced map $B\to C_*$ is adic.
Then $C_*$ is a perfectoid $A$-algebra.
\end{proposition}
\begin{proof} Let $I\subset B$ be any ideal of definition,
and $I_\bE\subset\bE(B)$ the corresponding ideal such that
$\bar u_A$ induces an isomorphism $\bE(B)/I_\bE\isom B/I$
(see remark \ref{rem_nice-topology}(iii)); then
$f^a\otimes_BB/I:\bE(B)/I_\bE\to C/IC$ is an \'etale morphism
of $(\bE,\fm_\bE)^a$-algebras. By \cite[Th.5.3.27]{Ga-Ra}, the
latter lifts uniquely to an \'etale morphism $g:\bE(B)^a\to C'$
of $(\bE,\fm_\bE)^a$-algebras such that $C'$ is almost finitely
presented as an $\bE(B)^a$-module. Since $\bE(B)$ is perfect,
\cite[Th.3.5.13]{Ga-Ra} implies that the Frobenius morphism
$$
(\Phi_\bE,\Phi_{C'}):((\bE(B),\fm_\bE),C')\to((\bE(B),\fm_\bE),C')
$$
is cartesian in the fibred category $\cB^a\Alg$ of
\cite[\S3.5.3]{Ga-Ra}. Since the functor $(-)_*$ on
$\cB^a$-algebras is cartesian (see \cite[\S3.5.4]{Ga-Ra}),
it follows that the morphism
$$
(\Phi_\bE,\Phi_{C'_*}):
((\bE(B),\fm_\bE),C'_*)\to((\bE(B),\fm_\bE),C'_*)
$$
is cartesian as well, {\em i.e.} $C'_*$ is a perfect
$\F_p$-algebra. Let us endow $C'_*$ with its $I_\bE$-adic
topology. We remark :

\begin{claim}\label{cl_etale-are-perf}
$C'_*$ is perfectoid.
\end{claim}
\begin{pfclaim} Since $\bE(B)^a$ is complete and separated
for its $I_\bE$-adic topology, \cite[Claim 5.3.25]{Ga-Ra},
implies that the natural morphism $C'\to\lim_{n\in\N}C'/I_\bE^nC'$
is an isomorphism, so the same holds for the induced map
$C'_*\to\lim_{n\in\N}(C'/I_\bE^nC')_*$. Especially, $C'_*$ is
complete and separated for the linear topology defined
by the system of ideals
$(J_n:=\Ker\,(C'_*\to(C'/I_\bE^nC')_*~|~n\in\N)$ (corollary
\ref{cor_limits-and-complete}(i)). Clearly
$I^n_\bE C'_*\subset J_n$ for every $n\in\N$, so $C'_*$ is also
complete and separated for its $I_\bE$-adic topology (lemma
\ref{lem_fontaine}), whence the claim.
\end{pfclaim}

By claim \ref{cl_etale-are-perf} it follows that there
exists a perfectoid $A$-algebra $D$ with an isomorphism
of perfectoid $\bE$-algebras $C'_*\isom\bE(D)$. We consider
the composition
$$
h:\bE(B)\xrightarrow{\ \eta_{\bE(B)}\ }
\bE(B)^a_*\xrightarrow{\ g_*\ }C'_*\isom\bE(D)
$$
where $\eta_{\bE(B)}$ is the unit of adjunction; since $h$
is an adic ring homomorphism, there exists a unique adic
map of $A$-algebras $k:B\to D$ with $\bE(k)=h$ (theorem
\ref{th_adic-to-adic}(i)). By construction, we have also
an isomorphism $(C/IC)^a\isom(D/ID)^a$ of $(B/I)^a$-algebras;
especially, $(k\otimes_BB/I)^a$ is a flat ring homomorphism,
and therefore the same holds for
$(k\otimes_BB/I^n)^a:(B/I^n)^a\to(D/I^nD)^a$, for every
$n\in\N$ (lemma \ref{lem_depressed}). By virtue of
\cite[Cor.3.2.11(ii,iii)]{Ga-Ra} and remark
\ref{rem_instead-of-3.2.1}, it follows that there exists
a unique isomorphism of $(B/I^n)^a$-algebras
$$
\omega_n:C/I^nC\isom(D/I^nD)^a
\qquad
\text{for every $n\in\N$}.
$$
Notice that the uniqueness of the morphisms $\omega_n$
implies that $\omega_{n+1}\otimes_BB/I^n=\omega_n$ for
every $n\in\N$. However, $D$ is complete and separated
for its $I$-adic topology, since $h$ is adic; the same
holds for $C$, since the latter is an almost projective
$B$-module of almost finite presentation
(\cite[Claim 5.3.25]{Ga-Ra}). Hence, the limit of the
inverse system $(\omega_n~|~n\in\N)$ is an isomorphism
$\omega:C\isom D^a$ of $B^a$-algebras, whence the
isomorphism $\omega_*:C_*\isom D^a_*$ of $B$-algebras,
and it suffices to invoke corollary \ref{cor_depressed}
to conclude.
\end{proof}

\sset\subsubsection{}\label{subsec_almost-purity}
Let now $\underline A:=(A,A^+,U)$ be any perfectoid
quasi-affinoid ring, and set
$$
\underline U:=\sSpec\,\underline A
\qquad
\underline U_\bE:=\bE(\underline U)
\qquad
X_A:=\Spec\,A
\qquad
Z_A:=X_A\setminus U
$$
(notation of \eqref{subsec_upgrade-E}). We let $\fm\subset A$
be the unique radical ideal such that $\Spec\,A/\fm=Z_A$. Also,
set $\bE:=\bE(A)$, $X_\bE:=\Spec\,\bE$, and let
$\fm_\bE\subset\bE$ be the unique ideal such that $\bar u_A$
induces a ring isomorphism $\bE/\fm_\bE\isom A/\fm$ (see remark
\ref{rem_nice-topology}(iii)). According to proposition
\ref{prop_down-to-business}(i), both pairs $(A,\fm)$ and
$(\bE,\fm_\bE)$ are basic setups.

\begin{theorem}\label{th_a-purity-for-perfectoids}
With the notation of \eqref{subsec_almost-purity}, the
pair $(X_A,Z_A)$ is normal and almost pure, relative to
the basic setup $(A,\fm)$.
\end{theorem}
\begin{proof} Let $j_A:U\to X_A$ be the open immersion;
since $j_{A*}\cO_U$ is a quasi-coherent $\cO_{\!X_A}$-algebra
(\cite[Ch.I, Prop.9.4.2(i)]{EGAI}), in order to check that
the pair $(X_A,Z_A)$ is normal, it suffices to show that the
induced map $A\to\cO_{\!X_A}(U)$ is almost injective and the
image of $A^a$ is integrally closed in $\cO_{\!X_A}(U)^a$.
These assertions follow from claim \ref{cl_perf-and-int-closed}(ii).

Next, pick any ideal of definition $I\subset A$,
set $\bE:=\bE(A)$, let $I_\bE\subset\bE$ be the unique
ideal of definition such that $\bar u_A$ induces a ring
isomorphism $\bE/I_\bE\isom A/I$, and consider the following
diagram of categories :
$$
\cD_A
\qquad :\qquad
{\diagram \bCov(U) & (A,\fm)^a\Et_\mathrm{afr}
\ar[l]_-{j_A^*} \ar[r]^-{i_A^*} &
(A/I,\fm/I)^a\Et_\mathrm{afr} \ddouble \\
\bCov(U_\bE) \ar@{-->}[u]^{w_A} &
(\bE,\fm_\bE)^a\Et_\mathrm{afr}
\ar[r]^-{i^*_\bE} \ar[l]_-{j^*_\bE} &
(\bE/I_\bE,\fm_\bE/I_\bE)^a\Et_\mathrm{afr}
\enddiagram}
$$
where $j_A^*$ (resp. $i^*_A$) is induced by the open
(resp. closed) immersion $j_A:U\to X_A$ (resp.
$i_A:Z_A\to X_A$) and likewise for $j^*_\bE$ and $i^*_\bE$.
We need to show that $j_A^*$ is an equivalence. However,
according to \cite[Th.5.5.7(iii)]{Ga-Ra}, both $i^*_A$
and $i^*_\bE$ are equivalences, and the same holds for
$j^*_\bE$, by virtue of theorem \ref{th_perfect-purity}(i).
Then, choose an arbitrary quasi-inverse functor $s_\bE$
for $j^*_\bE$ and $t_A$ for $i^*_A$; we let the dotted
arrow be the functor $w_A:=j^*_A\circ t_A\circ i^*_\bE\circ s_\bE:
\bCov(U_\bE)\to\bCov(U)$, so that $\cD_A$ is essentially
commutative, and we are reduced to checking that $w_A$
is an equivalence.

\begin{claim} We may assume that
$\underline A=\sGamma^\circ(\underline U)$
(notation of \eqref{subsec_Gamma-circ}).
\end{claim}
\begin{pfclaim} Say that
$\sGamma^\circ(\underline U)=(A_U^\circ,A^+_U,U)$,
set $X^\circ_U:=\Spec\,A^\circ_U$, and denote by
$\fm_U\subset A^\circ_U$ the unique radical ideal such that
$\Spec\,A^\circ_U/\fm_U=X^\circ_U\setminus U$.
Likewise, say that
$\sGamma^\circ(\underline U_\bE)=(\bE_U^\circ,\bE^+_U,U_\bE)$,
set $X^\circ_\bE:=\Spec\,\bE_U^\circ$ and denote by
$\fm_{\bE,U}\subset\bE^\circ_U$ the unique radical ideal such
that $\Spec\,\bE^\circ_U/\fm_{\bE,U}=X^\circ_\bE\setminus U_\bE$.
It follows easily that
$$
\fm_U\subset\fm A^\circ_U
\qquad
\fm_{U,\bE}\subset\fm_\bE\bE^\circ_U.
$$
Moreover, set $I_U:=IA^\circ_U$ and $I_{\bE,U}:=I_\bE\bE^\circ_U$;
the map $\bar u_{A_U^\circ}:\bE^\circ_U\to A^\circ_U$ induces a
ring isomorphism $\bE^\circ_U/I_{\bE,U}\isom A^\circ_U/I_U$.
Also, as in the foregoing, the projection
$A^\circ_U\to A^\circ_U/I_U$ and the open immersion
$j_{\bE,U}:U_\bE\to X^\circ_\bE$ induce equivalences of
categories
$$
i^*_U:(A^\circ_U,\fm_U)^a\Et_\mathrm{afr}\isom
(A^\circ_U/I_U,\fm_U/I_U)^a\Et_\mathrm{afr}
\qquad
j^*_{\bE,U}:(\bE^\circ_U,\fm_{\bE,U})^a\Et_\mathrm{afr}\isom
\bCov(U_\bE).
$$
After choosing quasi-inverse functors $t_U$ for $i^*_U$ and
$s_{\bE,U}$ for $j^*_{\bE,U}$, there follows an essentially
commutative diagram of functors
$$
\xymatrix{
\ar[d]^-{j^*_A} (A,\fm)^a\Et_\mathrm{afr} \ar@/_3pc/[dd] &
\ar[l]_-{t_A} (A/I,\fm/I)^a\Et_\mathrm{afr} \ar[dd] &
\ar[l]_-{i^*_\bE} (\bE,\fm_\bE)^a\Et_\mathrm{afr} \ar@/^3pc/[dd] \\
\bCov(U) & &
\bCov(U_\bE) \ar[u]^-{s_\bE} \ar[d]_-{s_{\bE,U}}  \\
\ar[u]_-{j^*_U} (A^\circ_U,\fm_U)^a\Et_\mathrm{afr} &
\ar[l]_-{t_U} (A^\circ_U/I_U,\fm_U/I_U)^a\Et_\mathrm{afr} &
\ar[l]_-{i^*_{\bE,U}}
(\bE^\circ_U,\fm_{\bE,U})^a\Et_\mathrm{afr}
}$$
whose unmarked vertical arrows are induced by the natural
ring homomorphisms $A\to A^\circ_U$, $A/I\to A^\circ_U/I_U$
and $\bE\to\bE^\circ_U$, and where $j^*_U$ is induced by
the open immersion $j_U:U\to X^\circ_U$. The claim follows
easily.
\end{pfclaim}

The restriction to the site $\cR(\underline U)$ of the stack
$\bCov^\wedge_{\underline U}$ of corollary \ref{cor_Cov-is-a-stack}
along the inclusion $\cR(\underline U)\subset\cQ(\underline U)$
is a stack that we denote again
$$
\pi_{\underline U}:\bCov^\wedge_{\underline U}\to\cR(\underline U).
$$
Recall the construction : for every rational subset
$R\subset\Spa\,\underline U$ we choose a complete and
separated quasi-affinoid scheme
$\underline U^\wedge_R:=(U^\wedge_R,\cT^\wedge_R,A^{\wedge+}_R)$
representing the sub-presheaf $h''_R$ of $h''_{\underline U}$;
then $\pi_{\underline U}^{-1}(R):=\bCov(U^\wedge_R)$. 
We shall consider similarly the quasi-affinoid schemes
$$
\underline X_A:=\sSpec\,(A,A^+,X_A)
\qquad\text{and}\qquad
\underline Y_A:=\sSpec\,A/I\otimes_A(A,A^+,X_A)
$$
(notation of example \ref{ex_shouldbe-quasi-affinoid}(i))
and the fibrations
$$
\pi_{\underline X_A^a}:
\bCov^\wedge_{(\underline X_A,\fm)^a}\to\cR(\underline U)
\qquad
\pi_{\underline Y_A^a}:
\bCov^\wedge_{(\underline Y_A,\fm)^a}\to\cR(\underline U)
$$
defined as follows. For any $R\in\Ob(\cR(\underline U))$,
we know that $(A^{\wedge\circ}_R,A^{\wedge+}_R,U^\wedge_R):=
\sGamma^\circ(\underline U^\wedge_R)$ is a perfectoid quasi-affinoid
ring (proposition \ref{prop_Letta-non-caduto}(iii) and
\eqref{subsec_Gamma-circ}), and we set
$X^\circ_R:=\Spec\,A^{\wedge\circ}_R$; then, let
$\fm_R\subset A^{\wedge\circ}_R$ be the unique radical ideal
such that $X^\circ_R\setminus U^\wedge_R=
\Spec\,A^{\wedge\circ}_R/\fm_R$. Set also $I_R:=IA^{\wedge\circ}_R$;
by proposition \ref{prop_down-to-business}(i), the pair
$(A^{\wedge\circ}_R,\fm_R)$ is a basic setup, and we let
$$
\pi^{-1}_{\underline X_A^a}(R):=(A^{\wedge\circ}_R,\fm_R)^a\Et_\mathrm{afp}
\qquad
\pi^{-1}_{\underline Y_A^a}(R):=
(A^{\wedge\circ}_R/I_R,\fm_R/I_R)^a\Et_\mathrm{afp}.
$$
If $R'\subset R$ is another rational subset, the
induced inclusion of sub-presheaves $h''_{R'}\subset h''_R$
corresponds to a morphism
$\underline U^\wedge_{R'}\to\underline U^\wedge_R$ of
quasi-affinoid schemes, which in turn yields a
morphism $\sGamma^\circ(\underline U^\wedge_R)\to
\sGamma^\circ(\underline U^\wedge_{R'})$ of perfectoid
quasi-affinoid rings, and it follows easily that
$$
\fm_{R'}\subset\fm_RA^{\wedge\circ}_{R'}
$$
whence, finally, induced functors
$$
j^*_{R'R}:\pi^{-1}_{\underline X_A^a}(R)\to\pi^{-1}_{\underline X_A^a}(R')
\qquad
j'^*_{R'R}:\pi^{-1}_{\underline Y_A^a}(R)\to\pi^{-1}_{\underline Y_A^a}(R').
$$
Moreover, for every further inclusion of rational
subset $R''\subset R'$ we have a natural isomorphisms
of functors
$$
j^*_{R''R'}\circ j^*_{R'R}\isom j^*_{R''R}
\qquad
j'^*_{R''R'}\circ j'^*_{R'R}\isom j'^*_{R''R}
$$
and the systems of such functors and isomorphisms of
functors amount to well defined pseudo-functors
$\cR(\underline U)^o\to\bCat$, whence the sought
fibrations $\pi_{\underline X_A^a}$ and $\pi_{\underline Y_A^a}$.
Furthermore, the system of open immersions
$j_R:U^\wedge_R\to X^\circ_R$ and the projections
$A^{\wedge\circ}_R\to A^{\wedge\circ}_R/I_R$ induce
pseudo-natural transformations
$$
(A^{\wedge\circ}_R/I_R,\fm_R/I_R)^a\Et_\mathrm{afp}
\xleftarrow{i^*_R}
(A^{\wedge\circ}_R,\fm_R)^a\Et_\mathrm{afp}
\xrightarrow{j^*_R}\bCov(U^\wedge_R)
\qquad
\text{for every $R\in\Ob(\cR(\underline U))$}
$$
whence morphisms of fibrations
$$
\bCov^\wedge_{(\underline Y_A,\fm)^a}
\xleftarrow{\ i^*_{\underline U}\ }\bCov^\wedge_{(\underline X_A,\fm)^a}
\xrightarrow{\ j^*_{\underline U}\ }\bCov^\wedge_{\underline U}.
$$
Likewise, by repeating the foregoing with the quasi-affinoid
schemes
$$
\underline X_\bE:=\sSpec\,(\bE,\bE^+,X_\bE)
\qquad
\underline Y_\bE:=\sSpec\,\bE/I_\bE\otimes_\bE(\bE,\bE^+,X_\bE)
$$
we get the following fibrations and morphisms of fibrations
over the site $\cR(\underline U_\bE)$ :
\set\begin{equation}\label{eq_the-latter}
\bCov^\wedge_{(\underline Y_\bE,\fm_\bE)^a}
\xleftarrow{\ i^*_{\underline U_\bE}\ }\bCov^\wedge_{(\underline X_\bE,\fm_\bE)^a}
\xrightarrow{\ j^*_{\underline U_\bE}\ }\bCov^\wedge_{\underline U_\bE}.
\end{equation}
However, in view of the isomorphism of sites
$\bE:\cR(\underline U)\isom\cR(\underline U_\bE)$ of remark
\ref{rem_tilt-of-site}(i), we may also regard \eqref{eq_the-latter}
as fibrations and morphisms of fibrations over the site
$\cR(\underline U)$. Moreover, by combining remark
\ref{rem_tilt-of-site}(ii) and \eqref{eq_osiris}, we obtain
natural isomorphisms
\set\begin{equation}\label{eq_prostaferesi}
\underline Y_A\times_{\underline X_A}\underline U^\wedge_R\isom
\underline Y_\bE\times_{\underline X_\bE}\underline U^\wedge_{\bE(R)}
\qquad
\text{for every $R\in\Ob(\cR(\underline U))$}
\end{equation}
and a simple inspection shows that the resulting ring isomorphism
$A^\circ_R/I_R\isom A^\circ_{\bE(R)}/I_{\bE(R)}$ maps $\fm_R/I_R$
isomorphically onto $\fm_{\bE(R)}/I_{\bE(R)}$. Summing up, we get
a natural fibrewise equivalence of fibrations over
$\cR(\underline U)$ :
\set\begin{equation}\label{eq_lusitania}
\bCov^\wedge_{(\underline Y_\bE,\fm_\bE)^a}\isom\bCov^\wedge_{(\underline Y_A,\fm)^a}.
\end{equation}
Furthermore, both $i^*_{\underline U}$ and $i^*_{\underline U_\bE}$ are
fibrewise equivalences, by virtue of \cite[Th.5.5.7(iii)]{Ga-Ra},
and the same holds for $j^*_{\underline U_\bE}$, by theorem
\ref{th_perfect-purity}(i). We may then find morphisms of fibrations
$$
s_{\underline U_\bE}:
\bCov^\wedge_{\underline U_\bE}\to\bCov^\wedge_{(\underline X_\bE,\fm_\bE)^a}
\qquad\text{and}\qquad
t_{\underline U}:\bCov^\wedge_{(\underline Y_A,\fm)^a}\to
\bCov^\wedge_{(\underline X_A,\fm)^a}
$$
that are quasi-inverse functors for $j^*_{\underline U_\bE}$ and
$i^*_{\underline U}$ (corollary \ref{cor_fibrations}(i)). The
composition of these functors yields a morphism of stacks
$$
w_{\underline U}:
\bCov^\wedge_{\underline U_\bE}\xrightarrow{\ s_{\underline U_\bE}\ }
\bCov^\wedge_{(\underline X_\bE,\fm_\bE)^a}\xrightarrow{\ i^*_{\underline U_\bE}\ }
\bCov^\wedge_{(\underline Y_\bE,\fm_\bE)^a}\isom\bCov_{(\underline Y_A,\fm)^a}
\xrightarrow{\ t_{\underline U}\ }
\bCov^\wedge_{(\underline X_A,\fm)^a}\xrightarrow{\ j^*_{\underline U}\ }
\bCov^\wedge_{\underline U}
$$
such that
$(w_{\underline U})_{\Spa\,\underline U}:\bCov(U_\bE)\to\bCov(U)$
is isomorphic to $w_A$.

For every $x\in\Spa\,\underline U$, let
$x_\bE:=\Spa(u_{\underline X})(x)\in\Spa\,\underline U_\bE$; by
propositions \ref{prop_faith-separ-cover} and
\ref{prop_stalks-for-stacks}, we are then reduced to checking
that the induced functor on stalks
$$
w_{\underline U}(x):\bCov^\wedge_{\underline U_\bE}(x_\bE)\to
\bCov^\wedge_{\underline U}(x)
$$
is an equivalence, for every such $x$. To this aim, suppose
first that $x$ is analytic; in this case, the completed
residue field $\kappa(x)^\wedge$ is a Tate valued field (see
\eqref{subsec_valuation-on-stalks}), so that the analytic
locus of $\Spec\,\kappa(x)^{\wedge+}$ is $\{\eta_x\}$, where
$\eta_x$ is the generic point. We have a well defined
quasi-affinoid ring
$$
\underline\kappa(x)^{\wedge\circ}:=
(\kappa(x)^{\wedge\circ},\kappa(x)^{\wedge+},\{\eta_x\})
$$
and the natural map $A\to\kappa(x)^{\wedge\circ}$ yields a
morphism of quasi-affinoid rings
$\pi^\wedge_x:\underline A\to\underline\kappa(x)^{\wedge\circ}$.
Likewise we get a morphism of quasi-affinoid rings
$$
\pi^\wedge_{x_\bE}:\bE(\underline A)\to
\underline\kappa(x_\bE)^{\wedge\circ}:=
(\kappa(x_\bE)^{\wedge\circ},\kappa(x_\bE)^{\wedge+},\{\eta_{x_\bE}\}).
$$
Moreover, by virtue of corollary \ref{cor_vals-are-perfectoid},
both $\underline\kappa(x)^{\wedge\circ}$ and
$\underline\kappa(x_\bE)^{\wedge\circ}$ are perfectoid, and
there exists a unique isomorphism of quasi-affinoid rings
\set\begin{equation}\label{eq_cassock}
\omega:\bE(\underline\kappa(x)^{\wedge\circ})\isom
\underline\kappa(x_\bE)^{\wedge\circ}
\qquad\text{such that}\qquad
\omega\circ\bE(\pi^\wedge_x)=\pi^\wedge_{x_\bE}.
\end{equation}
We may then repeat the foregoing constructions with
$\underline A$ replaced by $\underline\kappa(x)^{\wedge\circ}$
and $\underline\kappa(x_\bE)^{\wedge\circ}$ : we get first the
quasi-affinoid schemes
$$
\underline U(x):=\sSpec\,\underline\kappa(x)^{\wedge\circ}
\qquad
\underline U_\bE(x_\bE):=
\sSpec\,\underline\kappa(x_\bE)^{\wedge\circ}.
$$
Then we set $\underline A(x):=
(\kappa(x)^{\wedge\circ},\kappa(x)^{\wedge+},\Spec\,\kappa(x)^{\wedge\circ})$,
and
$$
\underline X_A(x):=
\sSpec\,\underline A(x)
\qquad
\underline Y_A(x):=\sSpec\,(A/I\otimes_A\underline A(x))
$$
and define likewise $\underline X_\bE(x_\bE)$ and
$\underline Y_\bE(x_\bE)$. Also, let
$\fm_x\subset\kappa(x)^{\wedge\circ}$ and
$\fm_{\bE,x}\subset\kappa(x_\bE)^{\wedge\circ}$ be the maximal
ideals. To these quasi-affinoid schemes we attach as in the
foregoing fibrations over the sites $\cR(\underline U(x))$
and $\cR(\underline U_\bE(x_\bE))$.
However, notice that the points of $\Spa\,\underline U(x)$
are (the equivalence classes of) the continuous valuations
$v:\kappa(x)^\circ\to\Gamma_{\!v\circ}$ with $\Ker\,v=0$ and
$v(\kappa(x)^+)\subset\Gamma^+_{\!v\circ}$; it is easily seen
that there is only one such valuation, {\em i.e.} the one
corresponding to the point $x$ of $\Spa\,\underline U$.
Likewise, we have $\Spa\,\underline U_\bE=\{x_\bE\}$; hence,
these fibrations amount to categories, or equivalently,
they can be identified naturally with their stalks over
the unique point $x$ of $\Spa\,\underline U$ and the unique
point $x_\bE$ of $\Spa\,\underline U_\bE$. Then, a simple
inspection yields an essentially commutative diagram of
categories :
$$
\cD'_A
\qquad : \qquad
{\diagram \bCov^\wedge_{(\underline Y_A,\fm)^a}(x) \ar[d] &
\bCov^\wedge_{(\underline X_A,\fm)^a}(x) \ar[l]_{i^*_{\underline U}(x)}
\ar[r]^-{j^*_{\underline U}(x)} \ar[d] &
\bCov^\wedge_{\underline U}(x) \ar[d] \\
\bCov^\wedge_{(\underline Y_A(x),\fm_x)^a} &
\bCov^\wedge_{(\underline X_A(x),\fm_x)^a} \ar[l]_{i^*_{\underline U(x)}}
\ar[r]^-{j^*_{\underline U(x)}} & \bCov^\wedge_{\underline U(x)}
\enddiagram}$$
whose right-most vertical arrow is an equivalence, by
virtue of \eqref{subsec_Florensky}. Likewise we have
an essentially commutative diagram $\cD'_\bE$ for the
corresponding fibrations over the site
$\cR(\underline U_\bE)$; moreover, in light of
\eqref{eq_cassock} we have as well an essentially
commutative diagram
$$
\xymatrix{ \bCov^\wedge_{(\underline Y_\bE,\fm_\bE)^a}(x_\bE) \ar[r] \ar[d] &
\bCov^\wedge_{(\underline Y_A,\fm)^a}(x) \ar[d] \\
\bCov^\wedge_{(\underline Y_\bE(x_\bE),\fm_{\bE,x})^a} \ar[r] &
\bCov^\wedge_{(\underline Y_A(x),\fm_x)^a}
}$$
whose vertical arrows are the same as the left-most vertical
arrows of $\cD'_A$ and $\cD'_\bE$, and whose horizontal
arrows are induced by the equivalences \eqref{eq_lusitania}.
Summing up, we are reduced to checking the following :

\begin{claim}
The functor $j^*_{\underline U(x)}$ is an equivalence.
\end{claim}
\begin{pfclaim} By unwinding the definitions, we see that
$\bCov^\wedge_{\underline U(x)}$ is the category of \'etale
$\kappa(x)^\wedge$-algebras, $\bCov^\wedge_{(\underline X_A(x),\fm_x)^a}$
is the category of \'etale
$(\kappa(x)^{\wedge\circ},\fm_x)^a$-algebras of almost finite
rank, and $j^*_{\underline U(x)}$ is the natural restriction
functor. In light of lemma \ref{lem_pure-almost-crit}(ii)
and remark \ref{rem_almost-pairs}(i), it follows already that
$j^*_{\underline U(x)}$ is fully faithful. Next, let $E$ be any
finite separable field extension of $\kappa(x)^\wedge$, and
$E^\circ$ the integral closure of $\kappa(x)^{\wedge\circ}$ in
$E$; it suffices to check that $E^{\circ a}$ is an \'etale
$\kappa(x)^{\wedge\circ a}$-algebra of finite rank. However,
notice that $\kappa(x)^{\wedge\circ}$ is a valuation ring of
rank one with field of fractions $\kappa(x)^\wedge$; moreover,
it is perfectoid for the ring topology such that the
natural map $A\to\kappa(x)^{\wedge\circ}$ is adic (corollary
\ref{cor_vals-are-perfectoid}(i)), and therefore it is
deeply ramified, by \cite[Prop.6.6.6]{Ga-Ra} and lemma
\ref{lem_perfectoid}(iv). Let also $E^\mathrm{s}$ be any
separable closure of $E$, and $E^{\mathrm{s}\circ}$ the
integral closure of $E^\circ$ in $E^\mathrm{s}$; by
\cite[Prop.6.6.2 and Rem.6.1.12(iv)]{Ga-Ra} it follows
that $(\Omega_{E^{\mathrm{s}\circ}/\kappa(x)^{\wedge\circ}})^a=0$, and
combining with \cite[Th.6.3.23]{Ga-Ra} we deduce that
$(\Omega_{E^\circ/\kappa(x)^{\wedge\circ}})^a=0$ (cp. the proof
of \cite[Prop.6.6.2]{Ga-Ra}). By invoking again
\cite[Th.6.3.23]{Ga-Ra}, it follows that
$\cD_{E^\circ/\kappa(x)^{\wedge\circ}}=E^{\circ a}$. Then the assertion
follows from \cite[Prop.6.3.8 and Lemma 4.1.27]{Ga-Ra}.
\end{pfclaim}

Lastly, in case $x$ is non-analytic, set
$$
Y(x):=\Spec\,(A/I\otimes_A\cO^\wedge_{\Spa\,\underline U,x})
\qquad
Y_\bE(x_\bE):=
\Spec\,(\bE/I_\bE\otimes_\bE\cO^\wedge_{\Spa\,\underline U_\bE,x_\bE})
$$
According to \eqref{subsec_Florensky} we have natural
equivalences of categories
$$
\Psi:\bCov^\wedge_{\underline U}(x)\isom\bCov(Y(x))
\qquad
\Psi_\bE:\bCov^\wedge_{\underline U}(x_\bE)\isom\bCov(Y_\bE(x_\bE)).
$$
Recall the construction of these equivalences : an
object of $\bCov_{\underline U}(x)$ is a pair $(R,\phi)$
where $R$ is a rational subset of $\Spa\,\underline U$
containing $x$, and $\phi:V\to U^\wedge_R$ is a finite
\'etale morphism of schemes; to such an object, the
equivalence $\Psi$ assigns the finite \'etale
morphism $Y(x)\times_{U^\wedge_R}\phi$. Similarly we may
describe $\Psi_\bE$. We get likewise well defined functors
$$
\bCov^\wedge_{(X_A,\fm_A)^a}(x)\xrightarrow{\ \Psi'\ }\bCov(Y(x))
\xleftarrow{\ \Psi''\ }\bCov^\wedge_{(\underline Y_A,\fm_A)^a}(x).
$$
Namely, an object of $\bCov^\wedge_{(X_A,\fm_A)^a}(x)$ is a pair
$(R,f:(A^{\wedge\circ}_R,\fm_R)^a\to B^a)$ where $R$
is as in the foregoing and $f$ is an \'etale covering
of almost rings, and notice that
$\fm_R\cO^\wedge_{\Spa\,\underline U,x}=\cO^\wedge_{\Spa\,\underline U,x}$,
so $f_x:=A/I\otimes_A\cO^\wedge_{\Spa\,\underline U,x}
\otimes_{A^{\wedge\circ}_R}f$ is a homomorphism of (usual) rings;
then the functor $\Psi'$ assigns to the pair $(R,f)$ the
finite \'etale covering $\Spec\,f_x$ of $Y(x)$. Similarly
we define $\Psi''$, and it is easily seen that the resulting
diagram of categories :
$$
\cD''_A
\qquad : \qquad
{\diagram \bCov^\wedge_{(\underline Y_A,\fm)^a}(x) \ar[rd]_{\Psi''} &
\bCov^\wedge_{(\underline X_A,\fm)^a}(x) \ar[l]_{i^*_{\underline U}(x)}
\ar[r]^-{j^*_{\underline U}(x)} \ar[d]^{\Psi'} &
\bCov^\wedge_{\underline U}(x) \ar[ld]^\Psi \\
& \bCov(Y(x))
\enddiagram}$$
is essentially commutative. Correspondingly, we have
an essentially commutative diagram $\cD''_\bE$ for the
stalks over $x_\bE$, and notice that all the arrows of
$\cD''_\bE$ are equivalences (details left to the reader).
On the other hand, the system of isomorphisms
\eqref{eq_prostaferesi} induces an isomorphism of schemes
$$
Y_A(x)\isom Y_\bE(x_\bE)
$$
whence an essentially commutative diagram
$$
\xymatrix{
\bCov^\wedge_{(\underline Y_\bE,\fm_\bE)^a}(x_\bE) \ar[r]^-{\Psi''_\bE} \ar[d] &
\bCov(Y_\bE(x_\bE)) \ar[d] \\
\bCov^\wedge_{(\underline Y_A,\fm)^a}(x) \ar[r]^-{\Psi''} & \bCov(Y(x))
}$$
whose left vertical arrow is induced by the equivalence
\eqref{eq_lusitania}, and the right vertical arrow is
an equivalence as well. We conclude that $\Psi''$ is an
equivalence; then the same holds for $j^*_{\underline U}(x)$,
and finally also for $w_{\underline U}(x)$, as sought. The
proof is complete.
\end{proof}

\sset\subsubsection{Almost purity in the formal perfectoid case}
We show now how to extend theorem \ref{th_a-purity-for-perfectoids}
to the case of formal perfectoid rings. We begin with some
preliminary observations :

\begin{lemma}\label{lem_basic-setup-on-form-perf}
Let $A$ be an adic topological ring that admits a finitely
generated ideal $I$ of adic definition, $A^\wedge$ the completion
of $A$, and $\fm\subset A$ an open ideal. The following holds :
\begin{enumerate}
\item
$\fm$ is a radical ideal of $A$ if and only if $\fm A^\wedge$
is a radical ideal of $A^\wedge$.
\item
$\fm=\fm^2$ if and only if $\fm A^\wedge=\fm^2 A^\wedge$.
\item
$\fm$ fulfills condition {\em({\bf B})} of\/ \cite[\S2.1.6]{Ga-Ra}
if and only if the same holds for $\fm A^\wedge$.
\end{enumerate}
\end{lemma}
\begin{proof} By assumption, $I^n\subset\fm$ for some
integer $n>0$; after replacing $I$ by $I^n$, we may then
assume that $I\subset\fm$.

(i): Let $\fn$ be the radical ideal of $\fm$; since
$I\subset\fm$, it is easily seen that $\fn/I$ is the
radical ideal of $\fm/I$ in the quotient $A/I$.
Since the completion map $A\to A^\wedge$ induces an
isomorphism $A/I\isom A^\wedge/IA^\wedge$, it follows that
$\fn A^\wedge/IA^\wedge$ is the radical ideal of
$\fm A^\wedge/IA^\wedge$; again, this easily implies that
$\fn A^\wedge$ is the radical ideal of $\fm A^\wedge$,
whence the assertion.

(ii): We have $\fm=\fm^2\Leftrightarrow\fm/I^2=\fm^2/I^2$,
and since the completion map induces an isomorphism
$A/I^2\isom A^\wedge/I^2A^\wedge$, we have
$\fm/I^2=\fm^2/I^2\Leftrightarrow
\fm A^\wedge/I^2A^\wedge=\fm^2/I^2A^\wedge$. Again, the
latter condition holds if and only if
$\fm A^\wedge=\fm^2 A^\wedge$, whence the assertion.

(iii): In view of (ii), it is clear that if $\fm$
fulfills condition ({\bf B}), then so does $\fm A^\wedge$.
Conversely, suppose that $\fm A^\wedge$ fulfills condition
({\bf B}); by (ii), it follows already that $\fm=\fm^2$,
and it remains to check that for every integer $k>0$, the
ideal $\fm$ is generated by the system
$x^k_\bullet:=(x^k~|~x\in\fm)$. However, let $a_1,\dots,a_n$
be a finite system of generators of $I$, and choose $N\in\N$
such that $I^N$ lies in the ideal generated by
$a^k_\bullet:=(a_1^k,\dots,a_n^k)$; by assumption $\fm A^\wedge$
is generated by the system $b^k_\bullet:=(b^k~|~b\in\fm A^\wedge)$,
hence the same holds for $\fm A^\wedge/I^NA^\wedge$.
The completion map induces an isomorphism
$\fm/I^N\isom\fm A^\wedge/I^NA^\wedge$ and maps $x_\bullet^k$
onto the family $b^k_\bullet$, so that $\fm/I^N$ is generated
by $x^k_\bullet$, and therefore the same holds for $\fm$, since
the family $x^k_\bullet$ contains $a^k_\bullet$.
\end{proof}

\begin{proposition}\label{prop_again-always-afr}
Let $A$ be a formal perfectoid ring, $I\subset A$ a
finitely generated ideal of adic definition,
$\fm\subset A$ a radical ideal with $I\subset\fm$, and
$U:=\Spec\,A\setminus\Spec\,A/\fm$. We have :
\begin{enumerate}
\item
$(A,\fm)$ is a basic setup fulfilling condition {\em({\bf B})}
(see \cite[\S2.1.1, \S2.1.6]{Ga-Ra}).
\item
For every subideal $J\subset I$, every almost finitely
generated almost projective $(A/J)^a$-module $P$ has almost
finite rank, for the almost structure induced by $(A,\fm)$.
\item
In the situation of {\em (ii)}, if moreover $U$ is quasi-compact,
then $P$ has finite rank.
\end{enumerate}
\end{proposition}
\begin{proof}(i) follows from lemma
\ref{lem_basic-setup-on-form-perf} and proposition
\ref{prop_down-to-business}(i).

(ii): After replacing $I$ by $I+pA$, we may assume that
$p\in I$. Set $X:=\Spec\,A/J$ and $U':=X\setminus\Spec\,A/I$.
Let $S\subset A$ be the multiplicative subset $1+I$, and
$B:=S^{-1}(A/J)$. Let also $\cP$ be the quasi-coherent
$\cO_X$-module associated with the $A/J$-module $P_!$. Then
$X=U'\cup\Spec\,B$, and $U'$ is quasi-compact, so that
$\cP_{|U'}$ is a locally free $\cO_{U'}$-module of finite
rank. Thus, it suffices to check that $P_B:=P\otimes_AB$
is almost projective $B^a$-module of almost finite rank.
Now, the ideal $I_B:=I/J$ of $B$ is obviously tight, and
by construction $I_B\subset\rad(B)$; moreover, the
localization $A/J\to B$ induces a ring isomorphism
$A/I\isom B/I_B$, whence an isomorphism of $(A/I)^a$-modules
$P_0:=P/IP\isom P_B/I_BP_B$. By lemma
\ref{lem_detect-afr-after-quot}(ii) we are then reduced to
checking that $P_0$ is an $(A/I)^a$-module of almost finite
rank. To this aim, let $\bE:=\bE(A)$, and pick a distinguished
element $\underline\alpha\in\Ker\,u_A$; by remark
\ref{rem_nice-topology}(ii), the map $\bar u_A$ induces
a ring isomorphism $\omega:\bE/\alpha_0\bE\isom A/pA$, and we
let $\fm_\bE$ and $I_\bE$ be the unique ideals of $\bE$ with
$$
\alpha_0\in\fm_\bE
\qquad
\alpha_0\in I
\qquad
\fm/pA\isom\fm_\bE/\alpha_0\bE
\qquad
I/pA\isom I_\bE/\alpha_0\bE.
$$
In particular, $I_\bE$ is finitely generated, and $(\bE,\fm_\bE)$
is a basic setup, by proposition \ref{prop_down-to-business}(i).
Via the isomorphism $\omega$, we may then regard $P_0$ as an
almost finitely generated almost projective $(\bE/I_\bE)^a$-module
(for the almost structure given by $(\bE,\fm_\bE)$), and clearly
it suffices to show that $P_0$ is of almost finite rank, when
regarded as an $(\bE/I_\bE)^a$-module. The latter follows from
proposition \ref{prop_always-afr}(i).

(iii): Arguing as in the proof of (ii), we reduce to showing
that if $U$ is quasi-compact, then $P_0$ is of finite rank,
when regarded as a $(\bE/I_\bE)^a$-module, and according to
proposition \ref{prop_always-afr}(ii), the latter will follow,
once we show that $Z:=\Spec\,\bE/\fm_\bE$ is constructible in
$\Spec\,\bE$. However,
$(\Spec\,A/pA)\setminus(\Spec\,A/\fm)=U\cap\Spec\,A/pA$ is
quasi-compact, hence the same holds for
$W:=(\Spec\,\bE/\alpha_0\bE)\setminus Z$. Thus,
$(\Spec\,\bE)\setminus Z=W\cup\Spec\bE[\alpha_0^{-1}]$ is
constructible in $\Spec\,\bE$, and then the same follows for $Z$.
\end{proof}

\sset\subsubsection{}\label{subsec_alm-pure-formal-perf}
Let $A$ be a formal perfectoid ring, $I\subset A$
a finitely generated ideal of adic definition, and
$\fm\subset A$ a radical ideal with $I\subset\fm$.
Set $X:=\Spec\,A$, $Z:=\Spec\,A/\fm$, and $U:=X\setminus Z$.
According to proposition \ref{prop_again-always-afr}(i),
we may consider the almost structure associated with the
basic setup $(A,\fm)$. We then have :

\begin{theorem}\label{th-alm-purity-form-perfectoid}
In the situation of \eqref{subsec_alm-pure-formal-perf},
let also $r\in\N$. The following holds :
\begin{enumerate}
\item
The pair $(X,Z)$ is almost pure.
\item
Let $\cM$ be any almost projective $\cO_{\!X}^a$-module of almost
finite rank. Then $\cM$ is a faithfully flat $\cO_{\!X}^a$-module
(resp. is an $\cO_{\!X}^a$-module of finite rank $\leq r$) 
if and only if the same holds for the $\cO_U$-module $\cM_{|U}$.
\end{enumerate}
\end{theorem}
\begin{proof}(i): Suppose first that $Z$ is constructible in
$X$. Let $A^\he$ and $A^\wedge$ be respectively the topological
henselization and the completion of $A$, and set
$X^\he:=\Spec\,A^\he$ and $X^\wedge:=\Spec\,A^\wedge$. Let
also $U^\he:=X^\he\times_XU$ and $U^\wedge:=X^\wedge\times_XU$.
In light of theorem \ref{th_a-purity-for-perfectoids}, it
suffices to check that both square subdiagrams of the following
induced diagram are $2$-cartesian :
$$
\xymatrix{ \cO^a_{\!X}\Et_\mathrm{afr} \ar[r] \ar[d] &
\cO^a_{\!X^\he}\Et_\mathrm{afr} \ar[d] \ar[r] &
\cO^a_{\!X^\wedge}\Et_\mathrm{afr} \ar[d] \\
\cO_{\!U}\Et_\mathrm{fr} \ar[r] &
\cO_{\!U^\he}\Et_\mathrm{fr} \ar[r] &
\cO_{\!U^\wedge}\Et_\mathrm{fr}.
}$$
Now, by \cite[Prop.5.5.6(i)]{Ga-Ra} it follows that the
natural functor
$$
\cO^a_{\!X}\Et_\mathrm{afr}\to\cO^a_{\!X^\he}\Et_\mathrm{afr}
\mathop{\times}^2_{\cO_{\!U^\he}\Et_\mathrm{fr}}\cO_{\!U}\Et_\mathrm{fr}
$$
is fully faithful. Thus, let $(\cE^\he,\cE_U,\xi)$ be the
datum of an object $\cE^\he$ of $\cO^a_{\!X^\he}\Et_\mathrm{afr}$,
an object $\cE_U$ of $\cO_{\!U}\Et_\mathrm{fr}$, and an
isomorphism $\xi$ in $\cO_{\!U^\he}\Et_\mathrm{fr}$ between
the images of $\cE^\he$ and $\cE_U$. By {\em loc.cit.}, there
exists a quasi-coherent $\cO^a_{\!X}$-algebra $\cE$ whose
restriction to $X^\he$ and to $U$ is isomorphic to $\cE^\he$
and respectively $\cE_U$, and in view of proposition
\ref{prop_again-always-afr}(iii), it remains to check that
$\cE$ is \'etale and almost finitely presented as an
$\cO^a_{\!X}$-module. Since the induced morphism
$U\sqcup X^\he\to X$ is quasi-compact and faithfully flat,
the assertion follows from lemmata \ref{lem_complem-qcoh}(i)
and \ref{lem_complem-qcoh-alg}(i). This shows that the left
square subdiagram of the foregoing diagram is $2$-cartesian.

Next, since both of the natural maps $A/I\to A^\he/IA^\he$
and $A/I\to A^\wedge/IA^\wedge$ are isomorphisms,
\cite[Th.5.5.7(iii)]{Ga-Ra} implies that the top horizontal
arrow of the right square subdiagram is an equivalence.
To conclude the proof in this case, it then suffices to show :

\begin{claim} The base change functor
$\cO_{U^\he}\Et_\mathrm{fr}\to\cO_{U^\wedge}\Et_\mathrm{fr}$
is an equivalence.
\end{claim}
\begin{pfclaim}
According to proposition \ref{prop_top-on-opens-fadic-case}(i)
there exists a unique f-adic topology on $B:=\cO_{U^\he}(U^\he)$
(resp. on $C:=\cO_{U^\wedge}(U^\wedge)$) such that the restriction
map $A^\he\to B$ (resp. $A^\wedge\to C$) is open. It follows
easily that the induced map $B\to C$ is continuous and adic,
since the same holds for the completion map $A\to A^\wedge$.
Set $A^{\he+}:=A^{\he\circ}\subset A^\he$ and
$A^{\wedge+}:=A^{\wedge\circ}\subset A^\wedge$; in view of lemma
\ref{lem_f-adics}(iii.a), the integral closure in $B$ (resp.
in $C$) of the image of $A^{\he+}$ (resp. of $A^{\wedge+}$) is a
subring of integral elements, that we denote $B^+$ (resp.
$C^+$). Recall also that the natural morphism of schemes
$\Spec\,B\to X^\he$ (resp. $X^\wedge\to\Spec\,C$) induces an
isomorphism $U_B:=U^\he\times_{X^\he}\Spec\,B\isom U^\he$ (resp.
$U_C:=U^\wedge\times_{X^\wedge}\Spec\,C\isom U^\wedge$); combining
with lemma \ref{lem_f-adic-pushout}, we deduce that the
resulting diagram
$$
\xymatrix{ \sSpec(C,C^+,U_C) \ar[r] \ar[d]_\phi &
\sSpec(B,B^+,U_B) \ar[d] \\
\sSpec(A^\wedge,A^{\wedge+},U^\wedge) \ar[r] &
\sSpec(A^\he,A^{\he+},U^\he)
}$$
is cartesian in the category of quasi-affinoid schemes.
Furthermore, the quasi-affinoid scheme $\sSpec(C,C^+,U_C)$
is complete, by virtue of lemma \ref{lem_topology-on-opens}(ii).
By arguing with universal properties, it follows easily
that the morphism $\phi$ identifies $\sSpec(C,C^+,U_C)$ with
the completion of $\sSpec(A^\wedge,A^{\wedge+},U^\wedge)$ : the
details are left to the reader. Then the claim follows from
theorem \ref{th_pass-to-completion} and proposition
\ref{prop_again-always-afr}(ii).
\end{pfclaim}

Lastly, we consider the case where $Z$ is an arbitrary closed
subset in $X$. Then, let $(I_\lambda~|~\lambda\in\Lambda)$ be
the filtered system of all finitely generated ideals of $A$
containing $I$ and contained in $\fm$. For every
$\lambda\in\Lambda$, let also $\fm_\lambda$ be the radical of
$I_\lambda$; then $Z_\lambda:=\Spec\,A/\fm_\lambda=\Spec\,A/I_\lambda$
and $U_\lambda:=X\setminus Z_\lambda$ are constructible subsets
of $X$, and $(A,\fm_\lambda)$ is a basic setup fulfilling
condition ({\em B}), for every such $\lambda$ (proposition
\ref{prop_again-always-afr}(i)). Both natural functors :
$$
(A,\fm)^a\Et_\mathrm{afr}\to
\Pslim{\lambda\in\Lambda}(A,\fm_\lambda)^a\Et_\mathrm{afr}
\qquad
\cO_{\!U}\Et_\mathrm{fr}\to
\Pslim{\lambda\in\Lambda}\cO_{\!U_\lambda}\Et_\mathrm{fr}
$$
are equivalences of categories (corollary
\ref{cor_limit-of-alm-struct}(ii,iv,v)), hence we are
reduced to prove the theorem for the pair $(X,Z_\lambda)$,
relative to the almost structure $(A,\fm_\lambda)$, for
every $\lambda\in\Lambda$. The latter is already known
by the previous case, so the proof is concluded.

(ii): Again, we assume first that $Z$ is constructible,
in which case $\cM$ has finite rank (proposition
\ref{prop_again-always-afr}(iii)). Let $A_\loc$ be the topological
localization of $A$ (see \eqref{subsec_localize-f-adic}); set
also $X_\loc:=\Spec\,A_\loc$. Then $A_\loc$ is a formal perfectoid
ring, and the localization $A\to A_\loc$ induces an isomorphism
of the respective completions $A^\wedge\isom A^\wedge_\loc$
(corollary \ref{cor_justify}(i)), and an isomorphism of
$A$-schemes $\pi:U_\loc:=X_\loc\times_XU\isom U$. Set $M:=\cM(X)$,
$M_\loc:=A_\loc\otimes_AM$, and let $\cM_\loc$ be the almost
projective $\cO^a_{\!X_\loc}$-module of finite rank corresponding
to $M_\loc$. Then $\pi$ induces a natural identification of
$\cM_{|U}$ with $(\cM_\loc)_{|U_\loc}$, so that the latter is a
faithfully flat $\cO_{U_\loc}$-module. We notice :

\begin{claim}\label{cl_a-purity-for-perfectoids}
In the situation of \eqref{subsec_almost-purity}, let $\cF$
be any almost projective $\cO_{\!X_A}^a$-module of finite rank.
Then $\cF$ is a faithfully flat $\cO_{\!X_A}^a$-module (resp.
is an $\cO_{\!X_A}^a$-module of finite rank $\leq r$) if and
only if the same holds for the $\cO_U$-module $\cF_{|U}$.
\end{claim}
\begin{pfclaim} Since the pair $(X,Z)$ is normal (theorem
\ref{th_a-purity-for-perfectoids}), the assertion follows
immediately from lemma \ref{lem_detect-faith-after-inj-ext}.
\end{pfclaim}

From claim \ref{cl_a-purity-for-perfectoids}, it follows
that $A^\wedge_\loc\otimes_{A_\loc}M_\loc$ is a faithfully flat
$(A^\wedge_\loc)^a$-module of finite rank; consequently
$M_\loc/IM_\loc$ is a faithfully flat $A^a_\loc/IA^a_\loc$-module
of finite rank. By construction, $I^a_\loc\subset\rad(A^a_\loc)$
is a tight ideal of $A^a_\loc$ (relative to the basic setup
$(A,\fm)$), hence $M_\loc$ is a faithfully flat $A^a_\loc$-module
(lemma \ref{lem_detect-afr-after-quot}(i)). Lastly, let $Y$
be the disjoint union of $U$ and $\Spec\,A_\loc$; we have an
obvious faithfully flat morphism of $A$-schemes $f:Y\to X$,
and the foregoing implies that $f^*\cM$ is a faithfully flat
$\cO^a_Y$-module, so $\cM$ is a faithfully flat
$\cO^a_{\!X}$-module, by lemma \ref{lem_complem-qcoh}(i).

Next, for a general $Z$, define the filtered system of
subideals $(\fm_\lambda~|~\lambda\in\Lambda)$ of $\fm$,
and the corresponding filtered system of open constructible
subsets $(U_\lambda~|~\lambda\in\Lambda)$ of $U$, as in the
proof of (i). Then if $\cM_{|U}$ is a faithfully flat
$\cO_U$-module (resp. if the rank of $\cM$ is $\leq r$),
obviously the same holds for the $\cO_{U_\lambda}$-module
$\cM_{|U_\lambda}$. For every $\lambda\in\Lambda$, denote by
$M_\lambda$ the image of $M$ in the category of
$(A,\fm_\lambda)^a$-modules; by the foregoing case, $M_\lambda$
is a faithfully flat $(A,\fm_\lambda)^a$-module (resp. has
rank $\leq r$) for every such $\lambda$. Then $M$ is a faithfully
flat $(A,\fm)^a$-module (resp. has rank $\leq r$), by corollary
\ref{cor_limit-of-alm-struct}(iii,iv).
\end{proof}

\subsection{Perfectoid Tate rings and perfectoid Abhyankar's lemma}
\label{sec_perf-Abhyankar}
In this section we introduce some special classes of
perfectoid quasi-affinoid rings that will play an important
role in the proof of Hochster's direct summand conjecture.

\begin{definition}\label{def_perf-Tate-ring}
Let $(A,\cT_A)$ be a Tate ring (see definition
\ref{def_bounded}(vi)).

(i)\ \
We say that $A$ is a {\em perfectoid Tate ring} (resp.
{\em formal perfectoid Tate ring}) if it has a subring of
definition $A_0$ that is perfectoid (resp. formal perfectoid),
for the topology induced by $\cT_A$ via the inclusion map $A_0\to A$.

(ii)\ \
We say that $A$ is {\em uniform}, if $A^\circ$ is a bounded
subset of $A$ (see definition \ref{def_bounded}(ii)). We let
$$
\Tate
\qquad\text{and}\qquad
\su.\Tate
$$
be the full subcategories of $\Z\tdu\TopAlg$ (see definition
\ref{def_top-ring}(iii)) whose objects are the Tate rings
and respectively the uniform Tate rings.
\end{definition}

\begin{lemma}\label{lem_uniform-Tate}
{\em (i)}\ \
Let $(A,\cT_A)$ be a Tate ring. The following conditions
are equivalent :
\begin{enumerate}
\alphaenu
\item
$A$ is uniform.
\item
The completion $A^\wedge$ of $A$ is uniform.
\item
$\cT_A$ agrees with the topology given by a
power-multiplicative semi-norm $\Vert\cdot\Vert:A\to\R$.
\end{enumerate}

{\em (ii)}\ \
Moreover, if the power-multiplicative semi-norm $\Vert\cdot\Vert$
defines the topology $\cT_A$, we have :
$$
A^\circ=\{a\in A~|~\Vert a\Vert\leq 1\}
\qquad\text{and}\qquad
A^{\circ\circ}=\{a\in A~|~\Vert a\Vert<1\}.
$$

{\em (iii)}\ \
The inclusion functor $\su.\Tate\to\Tate$\ \ admits a left
adjoint
$$
\su:\Tate\to\su.\Tate
\qquad
A\mapsto\su(A).
$$
We call $\su(A)$ the {\em uniformization} of the Tate ring $A$.
\end{lemma}
\begin{proof}(i.a)$\Rightarrow$(i.b): Indeed, if $A^\circ$
is a bounded subset of $A$, then its image $B$ in $A^\wedge$
is bounded in the latter ring; but then the same holds for
the topological closure $B^c$ of $B$ in $A^\wedge$ (remark
\ref{rem_something-on-bdd}(iii)). However, $B^c=(A^\wedge)^\circ$
(corollary \ref{cor_justify}(ii) and proposition
\ref{prop_replaces-Mat-Th.8.1}(iii)) whence the contention.

(i.b)$\Rightarrow$(i.a): Let $j:A\to A^\wedge$ be the completion
map; we have $j^{-1}((A^\wedge)^\circ)=A^\circ$, by corollaries
\ref{cor_justify}(ii) and \ref{cor_not-in-Bourbaki}(ii).
Since the topology of $A$ is induced by that of $A^\wedge$
via $j$, the assertion follows easily.

(i.a)$\Rightarrow$(i.c): Fix $\rho\in]0,1[$ and pick
$t\in A^\times\cap A^{\circ\circ}$. Notice that $A^\circ$ is
a subring of definition of $A$ (proposition
\ref{prop_f-adics}(ii)). Denote by
$\nu:A\to\Z\cup\{+\infty\}$ the order function associated
with the ideal $tA^\circ$ (see example \ref{ex_Samuel}(i)),
and set $|a|:=\rho^{\nu(a)}$ for every $a\in A$. Morever, let
$\Vert\cdot\Vert:A\to\R$ be the asymptotic Samuel function
of $tA^\circ$, {\em i.e.}
$\Vert a\Vert:=\displaystyle\lim_{n\to+\infty}|a^n|^{1/n}$ for
every $a\in A$. Notice that $\nu(t^ka)=k+\nu(a)$ for every
$k\in\N$ and $a\in A$, whence :
have :
\set\begin{equation}\label{eq_better-here}
\Vert t^ka\Vert=\rho^k\cdot\Vert a\Vert
\qquad\text{for every $k\in\N$ and $a\in A$}.
\end{equation}
By construction, $\Vert a\Vert\leq 1$ for every $a\in A^\circ$.
On the other hand, if $\Vert a\Vert>r>1$, then there exists
$k\in\N$ such that $|a^n|>r^n$ for every $n\geq k$, and
it follows easily that $a\notin A^\circ$. Combining with
\eqref{eq_better-here}, we conclude that
$t^kA^\circ=\{a\in A~|~\Vert a\Vert\leq\rho^k\}$ for every
$k\in\N$, whence (i.c).

(i.c)$\Rightarrow$(i.a): Clearly 
$\{a\in A~|~\Vert a\Vert<1\}=A^{\circ\circ}$. Pick
$t\in A^\times\cap A^{\circ\circ}$, so that
$r:=\Vert t\Vert\in]0,1[$; then $t^kA^{\circ\circ}$ is
an open ideal of $A^\circ$ for every $k\in\N$, and it lies
in the subset $\{a\in A~|~\Vert a\Vert<r^k\}$, whence (i.a).

(ii) follows by simple inspection of the definitions.

(iii): Let $A_0\subset A$ be a subring of definition,
$t\in A^\times\cap A^{\circ\circ}$, and $\nu:A\to\Z\cup\{+\infty\}$
the order function associated with the ideal $tA_0$.
For given $\rho\in]0,1[$, the map  $|\cdot|_A:A\to\R$
such that $|a|_A:=\rho^{\nu(a)}$ for every $a\in A$ is
a semi-norm, and the topology defined by $|\cdot|_A$
agrees with $\cT_A$. Then, let
$\Vert\cdot\Vert_A:A\to\R$ be the asymptotic Samuel
function of $tA_0$. Recall that $\Vert\cdot\Vert_A$ is a
multiplicative semi-norm, and let $\cT^\su_A$ be the topology
on $A$ defined $\Vert\cdot\Vert_A$. Lemma
\ref{lem_normalized-Samuel}(iii) implies that the identity
map of $A$ is a continuous ring homomorphism
$(A,\cT_A)\to(A,\cT^\su_A)$; also, the discussion of
example \ref{ex_Samuel} shows that $\cT^\su$ depends
only on the original topology $\cT_A$ (and is independent
of all auxiliary choices). Now, let $(B,\cT_B)$ be
any uniform Tate ring, and $f:(A,\cT_A)\to(B,\cT_B)$ a
continuous ring homomorphism. By lemma \ref{lem_deja-vu}(iv),
the map $f$ is f-adic, hence it restricts to an adic ring
homomorphism $f_0:A_0\to B^\circ$ (lemma
\ref{lem_f-adics}(iii.a,iii.b)), and the proof of (i) shows
that $\cT_B$ agrees with the topology defined by the
asymptotic Samuel function $\Vert\cdot\Vert_B:B\to\R$
of the ideal $f(t)\cdot B^\circ$. With this notation, it
follows easily that $f$ is a morphism of normed rings
$(A,|\cdot|_A)\to(B,\Vert\cdot\Vert_B)$ (see
\eqref{subsec_uniform-semi-normed}), and then lemma
\ref{lem_normalized-Samuel}(iv) shows that $f$ is
also a morphism of uniform semi-normed rings
$(A,\Vert\cdot\Vert_A)\to(B,\Vert\cdot\Vert_B)$.
In view of (i), we deduce that $f:(A,\cT^\su_A)\to(B,\cT_B)$
is a continuous ring homomorphism of uniform Tate rings.
Thus the sought functor is obtained by the rule :
$(A,\cT_A)\mapsto(A,\cT^\su_A)$.
\end{proof}

\begin{example}\label{ex_one-uniform-example}
Let $A$ be a Tate ring, $A_0\subset A$ a subring of definition,
and suppose there exists an integer $n>1$ such that
$\{x\in A~|~x^n\in A_0\}\subset A_0$. Then it follows easily
that $A^{\circ\circ}\subset A_0$, and since $A^{\circ\circ}$ is an
ideal of $A^\circ$, we deduce that $A$ is uniform.
\end{example}

\sset\subsubsection{}
\label{subsec_Perf-Tate-rings}
Let $(A,\cT_A)$ be a perfectoid Tate ring, and $A_0\subset A$
a perfectoid ring of definition, so that $A=A_0[t^{-1}]$, for
any $t\in A^\times\cap A_0^{\circ\circ}$; moreover, $tA_0$ is an
ideal of adic definition for $A_0$ (corollary \ref{cor_Tate}(iv)).
Let also $A^+$ be the integral closure of $A_0$ in $A$;
it follows that $U:=\Spec\,A$ is the analytic locus of
$\Spec\,A_0$, and
$$
\underline U:=(U,\cT_A,A^+)
$$
is a perfectoid affinoid scheme (see definition
\ref{def_perfectoid-qaff}(i)). Then, according to
\eqref{subsec_upgrade-E}, theorem
\ref{th_int-subrings-perfectoid}(ii) and corollary
\ref{cor_E-preserves-affines}(i), we get a perfectoid
affinoid scheme
$$
(U_\bE,\cT_\bE,\bE_U^+):=\bE(\underline U)
\qquad
\text{with $\bE_U:=\cO_{U_\bE}(U_\bE)$}
$$
so that $\bE_U$ is a complete f-adic ring, and its open
subring $\bE^+_U$ is perfectoid and naturally isomorphic
to $\bE(A^+)$. Let $(\alpha_n~|~n\in\N)$ be a distinguished
element in $\Ker\,u_{A^\circ}$; then the isomorphism of topological
rings $\bE^\circ_U/\alpha_0\bE^\circ_U\isom A^\circ/pA^\circ$ of
remark \ref{rem_nice-topology}(ii) induces an isomorphism
$$
\bar\omega:\bE^\circ_U/\bE^{\circ\circ}_U\isom A^\circ/A^{\circ\circ}.
$$
Recall also that the composition $\bE^\circ_U\to A$ of
$\bar u_{A^\circ}:\bE^\circ_U\to A^\circ$ with the localization
$A^\circ\to A$, factors through a continuous morphism of
monoids
$$
\phi^\flat_U:\bE_U\to A
$$
(lemma \ref{lem_phi-flat-map}(i) and proposition
\ref{prop_new-formula}(i)). Let
$\pi_A:A^\circ\to A^\circ/A^{\circ\circ}$ and
$\pi_\bE:\bE^\circ_U\to\bE^\circ_U/\bE^{\circ\circ}_U$ be the
projections; denote by $\cD_A$ and $\cD_\bE$ the set of
perfectoid subrings of adic definition of $A$ and respectively
$\bE_U$, and by $\cP_A$ the set of perfect subrings
of $A^\circ/A^{\circ\circ}$.

\begin{lemma}\label{lem_kobra}
In the situation of \eqref{subsec_Perf-Tate-rings},
the following holds :

{\em (i)}\ \
$(A,\cT_A)$ is a uniform Tate ring.

{\em (ii)}\ \
$(\bE_U,\cT_\bE)$ is a perfectoid Tate ring, and $U_\bE$
is the analytic locus of\/ $\Spec\,\bE^\circ_U$.

{\em (iii)}\ \
The rules : $B\mapsto\pi_A^{-1}B$ and
$B\mapsto(\bar\omega\circ\pi_\bE)^{-1}B$ establish natural
bijections
$$
\cD_\bE\stackrel{\sim}{\leftarrow}\cP_A\isom\cD_A.
$$
\end{lemma}
\begin{proof}(i) follows immediately from theorem
\ref{th_int-subrings-perfectoid}(i).

(ii): It is easily seen that $tA^\circ$ is an ideal of adic
definition for $A^\circ$, and $A^\circ[t^{-1}]=A$, hence $U$
is the analytic locus of $\Spec\,A^\circ$, and therefore
$U_\bE$ is the analytic locus of $\Spec\,\bE^\circ_U$. Pick
$m\in\N$ such that $p^m/t\in A^\circ$; by applying corollary
\ref{cor_Scholze-approx-lemma} to the quasi-affinoid
ring $(A^\circ,A^\circ,U)$, we deduce that there exists
$t'\in\bE^\circ_U$ such that
$$
v(t-\bar u_{A^\circ}(t'))\leq
v(p)\cdot\max(v(\bar u_{A^\circ}(t')),v(p^m))
\qquad
\text{for every $v\in\Cont(A^\circ)$}
$$
where $\bar u_{A^\circ}:\bE^\circ_U\to A^\circ$ is the natural
morphism of topological monoids as in \eqref{sec_Fontaine-only}.
It follows that $v(t)=v(\bar u_{A^\circ}(t'))$ for every
$v\in\Cont(A^\circ)$. We have $v(t)<1$ for every
$v\in\Cont(A^\circ)$ and $v(t)\neq 0$ for every $v\in\Cont(A)$,
since $t\in A^{\circ\circ}\cap A^\times$ (theorem
\ref{th_Cont-spectral}(i)); hence $w(t')<1$ for every
$w\in\Cont(\bE^\circ_U)$ (theorem \ref{th_Scholze-tilt}(i)),
and therefore $t'\in\bE^{\circ\circ}_U$ (corollary
\ref{cor_corcor}(ii)). Also, combining with
\eqref{subsec_what-name} we deduce that $w(t')\neq 0$
for every $w\in\Spa\,\bE(\underline U)$, hence
$t'\in\bE_U^\times$ (proposition \ref{prop_crit-invertible}).

(iii): For every $A'\in\cD_A$, the inclusion map
$A'\to A^\circ$ induces an open injective continuous
ring homomorphism $\bE':=\bE(A')\to\bE^\circ_U$, and
moreover we have $A^{\circ\circ}\subset A'$ and
$\bE^{\circ\circ}_U\subset\bE'$, by theorem
\ref{th_adic-to-adic}(i,iii.a). As in
\eqref{subsec_Perf-Tate-rings}, we deduce a ring isomorphism
$\bE'/\bE^{\circ\circ}_U\isom A'/A^{\circ\circ}$, which
shows that $A'/A^{\circ\circ}\in\cP_A$. Conversely, let
$B\in\cP_A$; according to proposition
\ref{prop_construct-new-perfs}, the open subring
$\pi_A^{-1}B=B\times_{A^\circ/A^{\circ\circ}}A^\circ$ of
$A^\circ$ is perfectoid for the topology induced by the
inclusion into $A^\circ$, and then $\pi^{-1}_AB\in\cD_A$.
Likewise one proves that the map $\cP_A\to\cD_\bE$ is
a bijection.
\end{proof}

\sset\subsubsection{}\label{subsec_Samuel-is-back}
Let now $(A,\cT_A)$ be a perfectoid Tate ring, $A_0\subset A$
a perfectoid subring of definition, $(B,\cT_B)$ an f-adic
ring, and $f:A\to B$ a continuous ring homomorphism.
Define the perfectoid affinoid schemes $\underline U$ and
$(U_\bE,\cT_\bE,\bE^+_U)$ and the f-adic ring $\bE_U$ as in
\eqref{subsec_Perf-Tate-rings}; by lemma \ref{lem_kobra}(ii) we
find $t'\in\bE_U^\times\cap\bE^{\circ\circ}_U$, and we let
$t:=\bar u_{A^\circ}(t')$, which is a topologically nilpotent
unit of $A$. Then more generally, the element
$t^\gamma:=\phi^\flat_U(t'^\gamma)\in A$ is well defined for every
$\gamma\in\Z[1/p]$. Notice that $s:=f(t)$ is a topologically
nilpotent unit of $B$, hence $B$ is a Tate ring; especially,
the non-analytic loci of both $\Spec\,A$ and $\Spec\,B$
are empty, so that $f$ is an f-adic ring homomorphism
(lemma \ref{lem_deja-vu}(iv)). By lemma
\ref{lem_f-adics}(iii.c) there exists a subring of definition
$B_0\subset B$ such that $f$ restricts to a continuous ring
homomorphism $f_0:A_0\to B_0$. Then $sB_0$ is an ideal of
adic definition of $B_0$, and $B=B_0[s^{-1}]$ (corollary
\ref{cor_Tate}(ii.b)). Let $s^\gamma:=f(t^\gamma)$ for every
$\gamma\in\Z[1/p]$. We fix $\rho\in]0,1[$ and consider the map
$$
|\cdot|_B:B\to\R\cup\{+\infty\}
\qquad
b\mapsto\inf\{\rho^\gamma~|~b\in s^\gamma B_0\}
$$
({\em cp.} \eqref{eq_back-on-track}). Let also
$\Vert\cdot\Vert_B:B\to\R\cup\{+\infty\}$ be the asymptotic
Samuel function of $sB_0$, as defined in example
\ref{ex_Samuel}(i). Recall furthermore that $A^\circ$ is
perfectoid for the topology induced by $A$ (theorem
\ref{th_int-subrings-perfectoid}(i)), and the pair
$(A^\circ,A^{\circ\circ})$ is a basic setup verifying condition
({\bf B}) of \cite[\S2.1.6]{Ga-Ra} (proposition
\ref{prop_down-to-business}(i)). Hence, the
$(A^\circ,A^{\circ\circ})^a$-algebra $B^a$ is well defined, and
in view of proposition \ref{prop_down-to-business}(ii) we see
that :
$$
B^a_*=\Hom_{A^\circ}(A^{\circ\circ},B)=
\Hom_A(A\otimes_{A^\circ}A^{\circ\circ},B)=\Hom_A(A,B)=B.
$$

\begin{lemma}\label{lem_Samuel-is-back}
{\em (i)}\ \
For every $b\in B$ we have :
$\displaystyle\lim_{n\to+\infty}|b^n|_B^{1/n}=\Vert b\Vert_B$.

{\em (ii)}\ \
$\Vert s^\gamma b\Vert_B=\rho^\gamma\cdot\Vert b\Vert_B$
for every $\gamma\in\Z[1/p]$ and every $b\in B$.

{\em (iii)}\ \
$B^{\circ\circ}=\{b\in B~|~\Vert b\Vert_B<1\}$.

{\em (iv)}\ \
$(B^\circ)^a_*=\{b\in B~|~\Vert b\Vert_B\leq 1\}$.
\end{lemma}
\begin{proof}(i): Recall that
$\Vert b\Vert_B:=\displaystyle\lim_{n\to+\infty}\rho^{\nu(b^n)/n}$,
where $\nu(b):=\inf\{k\in\Z~|~b\in s^kB_0\}$. For every $b\in B$,
set $\mu(b):=\sup\{\gamma\in\Z[1/p]~|~b\in s^\gamma B_0\}$, so
that $|b|_B=\rho^{\mu(b)}$. Clearly we have
$$
\nu(b^n)<\mu(b^n)<\nu(b^n)+1
\qquad
\text{for every $b\in B$ and $n\in\N$}.
$$
It follows that $\lim_{n\to+\infty}(|b^n|_B/\rho^{\nu(b^n)})^{1/n}=1$,
whence the contention.

(ii): It suffices to observe that
$|s^\gamma b|_B=\rho^\gamma\cdot|b|_B$ for every such $b$ and
$\gamma$, and apply (i).

(iii): It suffices to remark that
$\Vert b\Vert_B<1$ if and only if $x^n\in sB_0$ for some $n\in\N$.

(iv): Notice that $(B^\circ)^a_*\subset B^a_*=B$.
On the other hand, we have
$(B^\circ)^a_*=\Hom_{A^\circ}(A^{\circ\circ},B^\circ)$. It then follows
that $(B^\circ)^a_*=\{b\in B~|~bA^{\circ\circ}\subset B^\circ\}=
\{b\in B~|~bA^{\circ\circ}\subset B^{\circ\circ}\}$. Especially,
$s^\gamma b\in B^{\circ\circ}$ for every $b\in(B^\circ)^a_*$
and $\gamma\in\N[1/p]\setminus\{0\}$, whence
$\rho^\gamma\cdot\Vert b\Vert_B<1$ for every such $\gamma$,
in view of (ii) and (iii); {\em i.e.} $\Vert b\Vert_B\leq 1$.
Conversely, if $\Vert b\Vert_B\leq 1$, we get
$\Vert f(a)b\Vert_B\leq\Vert b\Vert_B\cdot\Vert f(a)\Vert_B<1$
for every $a\in A^{\circ\circ}$, due to (iii); then
$f(a)b\in B^{\circ\circ}$ for every such $a$, again due to (iii),
and finally, $b\in(B^\circ)^a_*$.
\end{proof}

\sset\subsubsection{}\label{subsec_pre-perfectoid}
Keep the situation of \eqref{subsec_Samuel-is-back}, and
set $C:=(B^\circ)^a_*$; lemma \ref{lem_Samuel-is-back}(ii,iv)
implies that
$$
s^\gamma C=\{b\in B~|~\Vert b\Vert_B\leq\rho^\gamma\}
\qquad
\text{for every $\gamma\in\Z[1/p]$}.
$$
In light of lemma \ref{lem_uniform-Tate}(iii), we may
now conclude that $C$ is the subring of power bounded
elements of the uniformization $\su(B)$ of $(B,\cT_B)$,
and the topology of $\su(B)$ is the f-adic topology such
that $C\subset B$ is a subring of definition, and $sC$
is an ideal of adic definition. Notice that
$f:(A,\cT_A)\to\su(B)$ is still an f-adic ring
homomorphism. Moreover, from lemmata
\ref{lem_uniform-Tate}(ii) and \ref{lem_Samuel-is-back}(iii)
we see that
\set\begin{equation}\label{eq_same-top-nil}
\su(B)^{\circ\circ}=B^{\circ\circ}.
\end{equation}

\begin{definition}
(i)\ \
With the notation of \eqref{subsec_Samuel-is-back}, we say
that $B$ is a {\em pre-perfectoid $A$-algebra}, if $\su(B)$
is a formal perfectoid Tate ring.
\end{definition}

\begin{remark}\label{rem_pre-perfectoid}
(i)\ \
With the notation of \eqref{subsec_pre-perfectoid}, endow
$C$ with its $sC$-adic topology $\cT_C$; from corollary
\ref{cor_justify}(ii) and theorem
\ref{th_int-subrings-perfectoid}(i) we see that $B$ is a
pre-perfectoid $A$-algebra if and only if $(C,\cT_C)^\wedge$
is a perfectoid ring.

(ii)\ \
Let $D\subset C$ be any subring that is $p$-integrally
closed in $C$ (see definition \ref{def_p-integrally-closed}(i)).
From \eqref{eq_same-top-nil} we see that :
$$
\text{$D$ is open in $(B,\cT_B)\Leftrightarrow
B^{\circ\circ}\subset D\Leftrightarrow
\su(B)^{\circ\circ}\subset D\Leftrightarrow$
$D$ is open in $\su(B)$}.
$$
Suppose then that $D$ is open in $(B,\cT_B)$, and denote
by $\cT_D$ the topology of $D$ induced by $\cT_C$; then we
claim that $B$ is a pre-perfectoid $A$-algebra if and only if
$(D,\cT_D)^\wedge$ perfectoid. Indeed, suppose that $B$ is a
pre-perfectoid $A$-algebra; since $B^{\circ\circ}\subset D$,
we see that $D/B^{\circ\circ}$ is a subring of $C/B^{\circ\circ}$,
and in view of \eqref{eq_same-top-nil} we get natural
identifications :
$$
C/B^{\circ\circ}=C/\su(B)^{\circ\circ}=\bar C:=
(C,\cT_C)^\wedge/(\su(B)^\wedge)^{\circ\circ}.
$$
Then, (i) implies that $C/B^{\circ\circ}$ is a perfect ring,
and $D/B^{\circ\circ}$ is a perfect subring, since $D$ is
$p$-integrally closed in $C$. In light of (i), lemma
\ref{lem_kobra}(iii) and proposition
\ref{prop_replaces-Mat-Th.8.1}(i,ii), we conclude that
$(D,\cT_D)^\wedge=D/B^{\circ\circ}\times_{\bar C}(C,\cT_C)$
is a perfectoid subring of definition of $\su(B)^\wedge$.
Conversely, suppose that $(D,\cT_D)^\wedge$ is perfectoid;
we know already that $D$ is open in $\su(B)$, and then
$(D,\cT_D)^\wedge$ is open in $\su(B)^\wedge$, and its
topology is induced by the inclusion into $\su(B)^\wedge$
(proposition \ref{prop_replaces-Mat-Th.8.1}(i,ii));
hence, $B$ is a pre-perfectoid $A$-algebra.

(iii)\ \
Lastly, let $D\subset B$ be an open subring, endow $D$
with the topology induced by $B$, and suppose that $D$
is formal perfectoid (definition \ref{def_formal-perfectoid});
then $D$ is $p$-integrally closed in $B$. Indeed, under
these assumptions, the topology of $D^\wedge$ is adic, so
$D^\wedge$ is a perfectoid subring of definition of $B^\wedge$,
hence $(B^\wedge)^{\circ\circ}\subset D^\wedge$, by virtue of lemma
\ref{lem_kobra}(iii); especially, $s^\gamma\in D$ for every
$\gamma\in\N[1/p]$, and $s^{-\gamma}p\in D$ for every sufficiently
small $\gamma\in\N[1/p]$ (corollary \ref{cor_not-in-Bourbaki}(ii)).
It then suffices to invoke the criterion of theorem
\ref{th_regular-seq-criterion}, together with lemma
\ref{lem_crit-p-integr-closed}(i).
\end{remark}

\begin{example}\label{ex_where-is-Manu}
Let us return to the situation of \eqref{subsec_alternative-meth} :
we consider a quasi-affinoid perfectoid ring
$\underline A:=(A,A^+,U_A)$ and a sequence $f_\bullet:=(f_0,\dots,f_n)$
of elements of $A_U:=\cO_{U_A}(U_A)$ that generates an open ideal
(for the f-adic topology on $A_U$ as in lemma
\ref{lem_topology-on-opens}(i)).

(i)\ \
Let $(\bE,\bE^+,U_\bE):=\bE(\underline A)$ (notation of
\eqref{subsec_upgrade-E}), and $\bE_U:=\cO_{U_\bE}(U_\bE)$.
As recalled in \eqref{subsec_Perf-Tate-rings}, the map
$\bar u_A:\bE\to A$ extends to a continuous morphism of
topological monoids $\phi^\flat_U:\bE_U\to A_U$; we may then find
$m\in\N$ and $e_0,\dots,e_n\in\bE_U$ verifying the inequalities
\eqref{eq_Naoko-1} and \eqref{eq_Naoko-2} for every $i=1,\dots,n$.
Let also $A_0\subset A_U$ be any subring of definition, and
$I_0\subset A_0$ an ideal of adic definition; recall that the
$I_0$-adic completion $B(f_\bullet)^\wedge$ of
$B(f_\bullet):=A_0[f_1/f_0,\dots,f_n/f_0]\subset A_U[1/f_0]$ is a
subring of definition of $\cO^\wedge_{\Spa\,\underline A}(S(f_\bullet))$,
with
$S(f_\bullet):=R_{A_U}\big(\frac{f_1}{f_0},\cdots,\frac{f_n}{f_0}\big)
\cap\Spa\,\underline A$. Moreover, the sequence
$\phi^\flat_U(e_0),\dots,\phi^\flat_U(e_n)$ generates
an open ideal of $A_U$ as well, according to proposition
\ref{prop_Scholze-tilt}(ii) and claim \ref{cl_extend-cor-taut-two},
and $S(f_\bullet)=R_{A_U}\big(\frac{\phi^\flat_U(e_1)}{\phi^\flat_U(e_0)},
\cdots,\frac{\phi^\flat_U(e_n)}{\phi^\flat_U(e_0)}\big)\cap
\Spa\,\underline A$, by proposition \ref{prop_Scholze-tilt}(iii).
It follows that the $I_0$-adic completion $B(e_\bullet)^\wedge$ of
$B(e_\bullet):=A_0\big[\frac{\phi^\flat_U(e_1)}{\phi^\flat_U(e_0)},
\dots,\frac{\phi^\flat_U(e_n)}{\phi^\flat_U(e_0)}\big]$ is another
subring of definition of $\cO^\wedge_{\Spa\,\underline A}(S(f_\bullet))$.
Therefore, the $p$-integral closures $D(f_\bullet)$ and $D(e_\bullet)$
of $B(f_\bullet)^\wedge$ and $B(e_\bullet)^\wedge$ in
$\cO^\wedge_{\Spa\,\underline A}(S(f_\bullet))$ both contain
$\cO^\wedge_{\Spa\,\underline A}(S(f_\bullet))^{\circ\circ}$.
Combining with proposition \ref{prop_Scholze-tilt}(iv), we
conclude that $D(f_\bullet)=D(e_\bullet)$.

(ii)\ \
In the situation of (i), notice that if
$f_0,\dots,f_n\in A$, we may take $e_0,\dots,e_n\in\bE$,
and by construction we shall then have
$f_i-\bar u_A(e_i)\in pA$ for $i=0,\dots,n$.
\end{example}

\begin{theorem}\label{th_crit-pre-perfectoid}
In the situation of \eqref{subsec_Samuel-is-back}, pick
$\gamma\in\N[1/p]\setminus\{0\}$ with $pB_0\subset s^{p\gamma}B_0$,
and suppose that the Frobenius endomorphism of
$B_0/s^{p\gamma}B_0$ is surjective. Then $B$ is pre-perfectoid.
\end{theorem}
\begin{proof} We exhibit an increasing chain of subrings
$$
B_0\subset B_1\subset B_2\subset\cdots
$$
of the $p$-integral closure of $B_0$ in $B$, verifying the
following conditions for every $i\in\N$ :
\begin{enumerate}
\alphaenu
\item
The Frobenius endomorphism $\Phi_i$ of $B_i/s^{p\gamma}B_i$
is surjective.
\item
Let $\bar\Phi_i:B_i/s^\gamma B_i\to B_i/s^{p\gamma}B_i$ be the
ring homomorphism induced by $\Phi_i$. Then the kernel of
the induced map $j_i:B_i/s^\gamma B_i\to B_{i+1}/s^\gamma B_{i+1}$
contains $\Ker\,\bar\Phi_i$.
\end{enumerate}
To this aim, we argue by induction on $i\in\N$. Hence, suppose
that for a given $i\in\N$ we have already exhibited a subring
$B_i$ fulfilling condition (a); we need to exhibit $B_{i+1}$
so that condition (b) holds as well for $B_i$, and moreover
such that condition (a) holds for $B_{i+1}$. Now, since (a) holds
for $B_i$, for every $x\in\Ker\,\bar\Phi_i$ we may find a sequence
$(a(x,n)~|~n\in\N)$ of elements of $B_i$, such that the
image of $a(x,0)$ in $B_i/s^\gamma B_i$ agrees with $x$, and :
$$
a(x,n+1)^p-a(x,n)\in s^{p\gamma}B_i
\qquad
\text{for every $n\in\N$}.
$$
It follows easily that $a(x,n)^{p^{n+1}}\in s^{p\gamma}B_i$, and
therefore $b(x,n):=a(x,n)\cdot s^{-\gamma/p^n}$ lies in the
$p$-integral closure of $B_i$ in $B$, for every $n\in\N$.
Let us then set
$$
B_{i+1}:=B_i[b(x,n)~|~(x,n)\in\Ker(\bar\Phi_i)\times\N]\subset B.
$$
Since the image of $b(x,0)\cdot s^\gamma$ in $B_i/s^\gamma B_i$
agrees with the image of $x$, we see that $j_i(x)=0$ for every
$x\in\Ker\,\bar\Phi_i$, hence condition (b) holds for $B_i$.
Lastly, notice that
$$
b(x,n+1)^p-b(x,n)\in B_i
\qquad
\text{for every $n\in\N$}.
$$
Since $B_i/s^{p\gamma}B_i\subset\Img\,\bar\Phi_{i+1}$, it follows
easily that condition (a) holds for $B_{i+1}$, as required.

Set $B_\infty:=\bigcup_{i\in\N}B_i=B_\infty$; by construction,
the Frobenius endomorphism of $B_\infty$ induces an isomorphism
$B_\infty/s^\gamma B_\infty\isom B_\infty/s^{p\gamma}B_\infty$, hence
$B_\infty$ is the $p$-integral closure of $B_0$ in $B$, by
lemma \ref{lem_crit-p-integr-closed}(i). Next, let
$B^\wedge_\infty$ be the completion of $B_\infty$ for its
$s$-adic topology; the Frobenius endomorphism of $B^\wedge_\infty$
induces an isomorphism $B^\wedge_\infty/s^\gamma B^\wedge_\infty\isom
B^\wedge_\infty/s^{p\gamma}B^\wedge_\infty$. Clearly $B^\wedge_\infty$
is a P-ring; by theorem \ref{th_regular-seq-criterion},
the ring $B_\infty^\wedge$ is then perfectoid, and the
assertion follows from remark \ref{rem_pre-perfectoid}(ii).
\end{proof}

\begin{corollary}
In the situation of \eqref{subsec_Samuel-is-back}, let
$I\subset B$ be any ideal, and endow $\bar B:=B/I$ with
the f-adic topology $\cT_{\bar B}$ induced by $\cT_B$ via
the projection $B\to\bar B$ (see example
{\em\ref{ex_f-adic-quotient}(iii)}). Then, if $B$ is a
pre-perfectoid $A$-algebra, the same holds for $\bar B$.
\end{corollary}
\begin{proof} Let $\cT^\su_B$ be the topology of $\su(B)$,
and denote by $\cT^\su_{\bar B}$ the f-adic topology induced
by $\cT^\su_B$ via the projection $B\to\bar B$. If $C$ is
any uniform Tate ring, and $f:(\bar B,\cT_B)\to C$ any
continuous ring homomorphism, it is easily seen that $f$
factors through the identity map
$\one_{\bar B}:(\bar B,\cT_B)\to(\bar B,\cT^\su_{\bar B})$, which
is continuous, and the continuous ring homomorphism
$f:(\bar B,\cT^\su_{\bar B})\to C$. By the universal property
of uniformization, this implies that $\one_{\bar B}$
induces an isomorphism of uniform Tate rings
$$
\su(\bar B,\cT_B)\isom\su(\bar B,\cT^\su_{\bar B})
$$
(whose underlying ring homomorphism is of course again
the identity map of $\bar B$). Thus, it suffices to check
that $(\bar B,\cT^\su_{\bar B})$ is pre-perfectoid. However,
since $B$ is pre-perfectoid, $C:=\su(B)^\circ$ is a perfectoid
ring, hence for every $\gamma\in\N[1/p]\setminus\{0\}$ such
that $pC\subset s^{p\gamma}C$, the Frobenius endomorphism of
$C/s^{p\gamma}C$ is surjective. By construction, the image
$\bar C$ of $C$ in $\bar B$ is a subring of definition
for the topology $\cT^\su_{\bar B}$, and it follows that
the Frobenius endomorphism of $\bar C/s^{p\gamma}\bar C$
is again surjective; then the assertion follows from
theorem \ref{th_crit-pre-perfectoid}.
\end{proof}

\sset\subsubsection{}\label{subsec_Zimbabwe}
Let $A$ be a formal perfectoid Tate ring, $A_0\subset A$ a
formal perfectoid subring of definition,
$t\in A_0^{\circ\circ}\cap A^\times$, and $P\in A_0[X]\setminus A_0$
a monic polynomial. We let
$B_0:=A_0[X^{1/p^\infty}]:=\bigcup_{n\in\N}A_0[X^{1/p^n}]$, and
$B:=A\otimes_{A_0}B_0$, and set
$$
C_{0,n}:=A_0[X^{1/p^n}]/(P)
\quad\text{for every $n\in\N$}\qquad
C_0:=B_0/PB_0
\qquad\text{and}\qquad
C:=B/PB.
$$
Clearly the natural map of $A_0$-algebras $C_{0,m}\to C_{0,n}$
is injective for every $n,m\in\N$ with $n\geq m$, and $C_{0,n}$
is a free $A_0$-module of rank $p^n\cdot\deg_X(P)$ for every
such $n$. Hence $C_0=\bigcup_{n\in\N}C_{0,n}$ is a faithfully flat
$A_0$-algebra, and the localization $C_0\to C$ is injective.
Lastly, let $D\subset C$  be the $p$-integral closure of $C_0$
in $C$ (see definition \ref{def_p-integrally-closed}(ii)), and
endow $D$ with its $t$-adic topology.

\begin{theorem}\label{th_Zimbabwe}
{\em (i)}\ \
In the situation of \eqref{subsec_Zimbabwe}, $D$ is a faithfully
flat $A_0$-algebra.

{\em (ii)}\ \
Moreover, $D$ is a formal perfectoid ring for its $t$-adic topology.
\end{theorem}
\begin{proof} Clearly $D[1/t]=C$ is a faithfully flat $A$-algebra,
and the image of $t$ is regular in $D$; by virtue of
\cite[Lemma 5.2.1]{Ga-Ra}, in order to show (i), it then suffices
to check that $D/tD$ is a faithfully flat $A_0/tA_0$-algebra. Recall
that the completion $A_0^\wedge$ is a subring of definition for the
completion $A^\wedge$ of $A$ (proposition
\ref{prop_complete-f-adic}(ii)); set $C'_0:=A^\wedge_0\otimes_{A_0}C_0$,
$C':=A^\wedge\otimes_AC$, and let $D'$ be the $p$-integral closure of
$C'_0$ in $C'$. Endow also $C_0$ (resp. $C'_0$) with its $t$-adic
topology, and $C$ (resp. $C'$) with the unique f-adic topology such
that $C_0$ (resp. $C'_0$) is a subring of definition; let $\cT_D$
(resp. $\cT_{D'}$) be the topology induced by $C$ on $D$ (resp. by
$C'$ on $D'$).
Taking into account proposition \ref{prop_complete-f-adic}, it
is easily seen that the natural maps $C_0\to C'_0$ and $C\to C'$
induce isomorphisms on the respective completions; by virtue of
lemma \ref{lem_crit-p-integr-closed}(ii), the same then holds for
the induced continuous ring homomorphism $(D,\cT_D)\to(D',\cT_{D'})$.

\begin{claim}\label{cl_Zimbabwe}
Let $D^\wedge$ be the completion of $(D,\cT_D)$. Then the
completion map $u:D\to D^\wedge$ induces an isomorphism
$D/t^kD\isom D^\wedge/t^kD^\wedge$ for every $k\in\N$.
\end{claim}
\begin{pfclaim} Let also $C_0^\wedge$ be the completion of $C_0$;
since $t^kC_0$ is an open $A_0$-submodule of $D$, the map $u$
induces an isomorphism of $A_0$-modules
$\bar u:D/t^kC_0\isom D^\wedge/t^kC^\wedge_0$. Hence,
$\bar u(t^kD/tC_0)=t^k\cdot\bar u(D/tC_0)=t^kD^\wedge/t^kC^\wedge_0$,
whence the assertion.
\end{pfclaim}

By applying claim \ref{cl_Zimbabwe} to both $D^\wedge$ and the
completion $D'^\wedge$ of $(D',\cT_{D'})$, we reduce to checking
that $D'/tD'$ is a faithfully flat $A_0/tA_0$-algebra, and that
$D'$ is perfectoid for its $t$-adic topology. Thus, we may replace
$A$ by $A^\wedge$, and assume that {\em $A$ is a perfectoid Tate ring}.

(ii): By construction, the Frobenius endomorphism of $C_0$
is a surjection; hence $C$ is a pre-perfectoid $A$-algebra,
by theorem \ref{th_crit-pre-perfectoid}. Then the assertion
follows from remark \ref{rem_pre-perfectoid}(ii).

(i): Define the affinoid schemes $\underline U:=(U,\cT_A,A^+)$
and $\underline U{}_\bE:=(U_\bE,\cT_\bE,\bE^+_U)$ as in
\eqref{subsec_Perf-Tate-rings}; recall that $\bE_0:=\bE(A_0)$
is a subring of definition of the Tate perfectoid ring
$\bE_U:=\cO_{U_\bE}(U_\bE)$ (lemma \ref{lem_kobra}(ii,iii)).
We may also assume that $t:=\bar u_{A_0}(t')\in A$ for some
$t'\in\bE_U^\times\cap\bE_0^{\circ\circ}$.

We endow $B_0$ with its $t$-adic topology,
and $B$ with the f-adic topology such that $B_0$ is a subring of
definition. We have then a short exact sequence of flat $A_0$-modules
$$
\Sigma \qquad : \qquad
0\to PB_0\to B_0\xrightarrow{\pi_0}C_0\to 0.
\qquad\qquad
$$
As in \eqref{subsec_Samuel-is-back}, we let
$t^\gamma:=\phi^\flat_U(t'^\gamma)\in A$ for every $\gamma\in\Z[1/p]$.
With this notation, the sequence $A_0/t^\gamma A_0\otimes_{A_0}\Sigma$
is still short exact, for every $\gamma\in\N[1/p]$, which means that
$t^\gamma B_0\cap PB_0=t^\gamma PB_0$ for every such $\gamma$.
Especially, the topology of $B_0$ induces the $t$-adic topology
on the ideal $PB_0$, and after taking completions, we still
get a short exact sequence :
$$
\Sigma^\wedge \qquad : \qquad
0\to PB_0^\wedge\to B^\wedge_0
\xrightarrow{\pi^\wedge_0}C_0^\wedge\to 0.
\qquad\qquad
$$
Recall that the completion $C^\wedge$ of $C$ is naturally identified
with $C_0^\wedge[t^{-1}]$, and $C_0^\wedge$ is a subring of definition
for $C^\wedge$ (proposition \ref{prop_complete-f-adic}(ii,iii)).
Likewise, $B_0^\wedge$ is a subring of definition of
$B^\wedge=B_0^\wedge[t^{-1}]$. We consider the affinoid ring
$\underline B^\wedge:=(B^\wedge,B^{\wedge\circ})$, and the rational
subset
$$
R_n:=
R_{B^\wedge}\Big(\frac{P}{t^n}\Big)\cap\Spa\,\underline B^\wedge
\qquad
\text{for every $n\in\N$}.
$$
Endow the subring $B_n:=B_0^\wedge[P/t^n]\subset B^\wedge$
with its $t$-adic topology, and let $\cT_n$ be the f-adic
topology on $B^\wedge$ such that $B_n$ is a subring of
definition; then the completion $B_n^\wedge$ of $B_n$ is a
subring of definition of
$\cO^\wedge_{\Spa\,\underline B^\wedge}(R_n)=(B^\wedge,\cT_n)^\wedge$,
for every $n\in\N$.
Notice that the projection $\pi_0^\wedge$ extends to a continuous
open ring homomorphism $\pi^\wedge:B^\wedge\to C^\wedge$, and for
every $n\in\N$, the image of the induced continuous map
$$
\Spa\,\pi^\wedge:
\Spa(C^\wedge,C^{\wedge\circ})\to\Spa\,\underline B^\wedge
$$
lies in $R_n$. It follows that $\pi^\wedge$ factors uniquely
through a continuous ring homomorphism
$$
\psi_n:\cO^\wedge_{\Spa\,\underline B^\wedge}(R_n)\to C^\wedge
\qquad
\text{for every $n\in\N$}
$$
(see remark \ref{rem_yoneda-rationals}(i)). Moreover, we
have $R_{n+1}\subset R_n$ for every $n\in\N$, and clearly
the composition of $\psi_{n+1}$ with the restriction map
$\cO^\wedge_{\Spa\,\underline B^\wedge}(R_n)\to
\cO^\wedge_{\Spa\,\underline B^\wedge}(R_{n+1})$ agrees with
$\psi_n$. Furthermore, a direct inspection shows that
$\psi_n(B^\wedge_n)=C_0^\wedge$ for every $n\in\N$.
Summing up, we deduce a commutative diagram of rings :
\set\begin{equation}\label{eq_danced-all-night}
{\diagram
B':=\displaystyle{\colim_{n\in\N}B^\wedge_n} \ar[r] \ar[d]_\alpha &
\cO:=\displaystyle{\colim_{n\in\N}\cO^\wedge_{\Spa\,\underline B^\wedge}(R_n)}
\ar[d]^\beta \\
C_0^\wedge \ar[r] & C^\wedge
\enddiagram}
\end{equation}
whose vertical (resp. horizontal) arrows are surjections
(resp. injections).

\begin{claim}\label{cl_evier-bouche}
The diagram \eqref{eq_danced-all-night} is cartesian.
\end{claim}
\begin{pfclaim} The assertion means that the induced map
$\Ker\,\alpha\to\Ker\,\beta$ is bijective. However, this
map is clearly injective, so it remains to check its
surjectivity. To this aim, notice that for every $n\in\N$
and every $g_n\in\cO_{\Spa\,\underline B^\wedge}(R_n)$ there exists
$N\in\N$ such that
$$
t^Ng_n=\sum_{i=0}^\infty b_i\cdot\Big(\frac{P}{t^n}\Big)^i
\qquad
\text{with $b_i\in B_0^\wedge$ for every $i\in\N$ and
$\lim_{i\to+\infty}b_i=0$}.
$$
With this notation, it follows easily that $g_n\in\Ker\,\psi_n$
if and only if $b_0=0$. Suppose then that $\psi_n(g_n)=0$; the
image $g_{n+N}$ of $g_n$ in $\cO_{\Spa\,\underline B^\wedge}(R_{n+N})$
can be written as :
$$
g_{n+N}=\sum_{i=1}^\infty t^{N(i-1)}b_i\cdot\Big(\frac{P}{t^{n+N}}\Big)^i
$$
and we have $t^{N(i-1)}b_i\in B^\wedge_0$ for every $i\geq 1$, and
$\displaystyle{\lim_{i\to+\infty}}t^{N(i-1)}b_i=0$. Hence
$g_{n+N}\in B^\wedge_{n+N}$, and clearly $g_{n+N}\in\Ker\,\psi_{n+N}$,
whence the contention.
\end{pfclaim}

Let $\cD$ be the $p$-integral closure of $B'$ in $\cO$,
and endow $\cD$ with its $t$-adic topology; it follows
easily from claim \ref{cl_evier-bouche} that
$\cD=\beta^{-1}D$. Moreover, the induced ring homomorphism
$\beta_{|\cD}:\cD\to D$ is surjective, and its kernel is
$\Ker\,\beta$. The latter is clearly $t$-divisible, hence
$\beta_{|\cD}$ induces an isomorphism of $A_0$-algebras
$$
\cD/t^n\cD\isom D/t^nD
\qquad
\text{for every $n\in\N$}.
$$
Thus, it suffices to check that the natural ring
homomorphism $A_0\to\cD$ is adically faithfully flat.
To this aim, let also $\cD_n$ be the $p$-integral closure
of $B^\wedge_n$ in $\cO^\wedge_{\Spa\,\underline B^\wedge}(R_n)$,
for every $n\in\N$, and endow $\cD_n$ with its $t$-adic
topology; in view of lemma
\ref{lem_crit-p-integr-closed}(iii.b) we are reduced to
checking that $\cD_n$ is an adically faithfully flat
$A_0$-algebra for every $n\in\N$.

Now, notice that $B_0^\wedge$ is a perfectoid ring, by virtue
of lemma \ref{lem_graded-perfectoid}, and we have a natural
identification of topological rings :
$$
\bE(B^\wedge_0)\isom\bE_0[X^{1/p^\infty}]^\wedge
\qquad\text{such that}\qquad
\bar u_{B_0^\wedge}(X^\gamma)=X^\gamma
\qquad\text{for every $\gamma\in\N[1/p]$}
$$
where $\bE_0[X^{1/p^\infty}]^\wedge$ denotes the $t'$-adic
completion of $\bE_0[X^{1/p^\infty}]$. Next, let
$(\alpha_n~|~n\in\N)$ be a distinguished element in the
kernel of $u_{A_0}:W(\bE_0)\to A_0$; after replacing
$t'$ by $t'^\gamma$ for some sufficiently small
$\gamma\in\N[1/p]\setminus\{0\}$, we may assume that
$\alpha_0\bE_0\subset t'\bE_0$, in which case the isomorphism
$\omega$ of remark \ref{rem_nice-topology}(ii) induces ring
isomorphisms
\set\begin{equation}\label{eq_usual-E_0-A_0-game}
\bE_0/t'\bE_0\isom A_0/tA_0
\qquad
\bE(B^\wedge_0)/t'\bE(B^\wedge_0)\isom B^\wedge_0/tB^\wedge_0.
\end{equation}
On the other hand, set $U_B:=\Spec\,B^\wedge$; then
$(B^\wedge_0,B^\wedge_0,U_B)$ is a perfectoid quasi-affinoid
ring, and by example \ref{ex_where-is-Manu} we may then
find $e_n\in\bE(B^\wedge_0)$ such that the $\cD_n$ is the
$p$-integral closure in $\cO^\wedge_{\Spa\,\underline B^\wedge}(R_n)$ of
the $t$-adic completion of $B^\wedge_0[\bar u_{B^\wedge_0}(e_n)/t^n]$,
and moreover -- after replacing again $t'$ by $t'^\gamma$ for
some small $\gamma\in\N[1/p]\setminus\{0\}$ -- we may assume
that $P-\bar u_{B^\wedge_0}(e_n)\in tB^\wedge_0$. Hence, there exist
a monic polynomial $Q_n\in\bE_0[X]$ and
$g_n\in\bE_0[X^{1/p^\infty}]^\wedge$ with
\set\begin{equation}\label{eq_present-e_n}
e_n=Q_n+t'g_n.
\end{equation}
According to proposition \ref{prop_Letta-non-caduto}(iii),
the subring
$$
D_n:=B_0^\wedge[\bar u_{B^\wedge_0}(e_n^{1/p^k})/t^{n/p^k}~|~k\in\N]
\subset B^\wedge
$$
is formal perfectoid for the topology $\cT'_n$ induced by
$\cT_n$. Then $\cT'_n$ is the $t$-adic topology, and combining
with remark \ref{rem_pre-perfectoid}(iii) and lemma
\ref{lem_crit-p-integr-closed}(ii), we deduce that
$$
(D_n,\cT'_n)^\wedge=\cD_n.
$$
Thus, we are further reduced to checking that the $A_0$-algebra
$D_{n,k}:=B^\wedge_0[\bar u_{B^\wedge_0}(e_n^{1/p^k})/t^{n/p^k}]$
is adically faithfully flat for every $n,k\in\N$ with $n>0$
and $n/p^k\leq 1$. Since $t$ is a regular element of $D_{n,k}$,
the local flatness criterion (\cite[Th.22.3]{Mat}) further
reduces to showing that $D_{n,k}/(t^{n/p^k})$ is a faithfully
flat $A_0/(t^{n/p^k})$-algebra, for every such $n$ and $k$. To
this aim, notice that \eqref{eq_usual-E_0-A_0-game} induces
an isomorphism $E_{n,k}:=\bE_0/(t'^{n/p^k})\isom A_0/(t^{n/p^k})$,
and $B^\wedge_0/(t^{n/p^k})$ is isomorphic to the $A_0$-algebra
$E_{n,k}[X^{1/p^\infty}]$ (as $n/p^k\leq 1$). Also, the image of
$\bar u_{B^\wedge_0}(e_n^{1/p^k})$ in $E_{n,k}[X^{1/p^\infty}]$ is the
regular element $Q_n^{1/p^k}$, due to \eqref{eq_present-e_n}.
Thus, the sequence $(t^{n/p^k},\bar u_{B^\wedge_0}(e_n^{1/p^k}))$ is
regular in the ring $B^\wedge_0$; in particular, this sequence
is completely secant (proposition \ref{prop_Kosz-cptl-sec}),
so the natural map
$$
B^\wedge_0[Y]/(\bar u_{B^\wedge_0}(e_n^{1/p^k})-t^{n/p^k}Y)\to D_{n,k}
$$
is an isomorphism of $A_0$-algebras (lemma \ref{lem_from-SGA6}).
We are then furter reduced to checking that for every $n,k,r\in\N$
with $n>0$, $r\geq k$ and $n/p^k\leq 1$, the $E_{n,k}$-algebra
$$
E_{n,k}[X^{1/p^r},Y]/(Q_n^{1/p^k})
$$
is faithfully flat. The latter follows from
\cite[Ch.IV, Prop.11.3.7]{EGAIV-3}.
\end{proof}

\sset\subsubsection{}\label{subsec_make-a-comp-pair-ofalm-str}
Let $A$ be a perfectoid Tate ring, $A_0\subset A$ a perfectoid
subring of definition, and $(A_0,\fm)$ a basic setup in the sense
of \cite[\S2.1.1]{Ga-Ra}, such that $\fm$ is the filtered union
of a countable system of principal subideals. If $\fm$ is an open
radical ideal of $A_0$, then we have attached to $\fm$ an open
radical ideal $\fm_\bE$ of $\bE_0:=\bE(A_0/pA_0)$, so that
$(\bE_0,\fm_\bE)$ is a basic setup as well (see proposition
\ref{prop_down-to-business}). One can show that here $\fm$
is still a radical ideal; however, it is not necessarily open,
but nevertheless we can construct a useful basic setup
$(\bE_0,\fm_\bE)$ as follows. By assumption, $\fm$ is the filtered
union of a countable system of principal ideals; hence, we may
find a generating system $(a_n~|~n\in\N)$ for $\fm$ such that for
every $n\in\N$ there exists an element $c_n\in A$ with
$a_n=c_na^p_{n+1}$. For every $x\in A_0$, let $\bar x\in A_0/pA_0$
be the class of $x$; since $A_0$ is perfectoid, for every $n\in\N$
there exists $(\bar x_{n,k}~|~k\in\N)\in\bE_0$ with
$\bar x_{n,0}=\bar c_n$ (recall that $\bar x{}_{n,k+1}^p=\bar x_{n,k}$
for every $k\in\N$). Set
$$
\alpha_{n,k}:=
\bar a_{n+k}\cdot\bar x_{n,k}\cdot\bar x_{n+1,k-1}\cdots\bar x_{n+k-1,1}
\qquad
\text{for every $n\in\N$}.
$$
Thus, $\alpha_{n,0}=\bar a_n$ and it is easily seen that
$\alpha_{n,\bullet}:=(\alpha_{n,k}~|~k\in\N)\in\bE_0$ for
every $n\in\N$; moreover,
$\alpha_{n,\bullet}\in\alpha^p_{n+1,\bullet}\bE_0$ for every $n\in\N$.
Let then $\fm_\bE\subset\bE_0$ be the ideal generated by
the system $(\alpha_{n,\bullet}~|~n\in\N)$; it follows that
$(\bE_0,\fm_\bE)$ is a basic setup.

\begin{remark}\label{rem_make-comp-pair-alm-struct}
(i)\ \
In the situation of \eqref{subsec_make-a-comp-pair-ofalm-str},
we have :
\set\begin{equation}\label{eq_same-p-top-closure}
\bigcap_{n\in\N}(\bar u_{A_0}(\fm_\bE)\cdot A_0+p^nA_0)=
\bigcap_{n\in\N}(\fm+p^nA_0).
\end{equation}
Indeed, on the one hand, $\bar u_{A_0}(\alpha_{n,\bullet})$ is the
limit of the $p$-adically convergent sequence
$$
((a_{n+k}x_{n,k}\cdots x_{n+k-1,1})^{p^k}~|~k\in\N)
$$
hence it lies in the $p$-adic topological closure of $\fm$.
On the other hand, by construction for every $x\in\fm$ there
exists $y\in\bar u_{A_0}(\fm_E)$ such that $x-y\in pA_0$; then
$x^{p^n}-y^{p^n}\in p^nA_0$ for every $n\in\N$ (lemma
\ref{lem_basic-cong}(i)), and since $\fm$ is generated by
$(x^{p^n}~|~x\in\fm)$, we deduce the converse inclusion.

(ii)\ \
Quite generally, let $(V,\fm)$ be a basic setup in the sense
of \cite[\S2.1.1]{Ga-Ra}, such that $\fm$ is a filtered countable
union of principal subideals; then it is easily seen that we may
write $\fm=\bigcup_{n\in\N}Va_n$ for a system $(a_n~|~n\in\N)$
of elements of $A$ such that $Va_n\subset\fm\cdot a_{n+1}$
for every $n\in\N$ : the details are left to the reader.
Hence, for every $n\in\N$ let $c_n\in\fm$ such that
$a_n=c_na_{n+1}$. We consider the inductive system
$L_\bullet:=(L_n~|~n\in\N)$ with $L_n:=V$ for every $n\in\N$,
with $V$-linear transition map $\phi_n:L_n\to L_{n+1}$ such
that $\phi_n(1):=c_n$ for every $n\in\N$. The inductive limit
$L$ of $L_\bullet$ is the set of equivalence classes $[x,n]$
of all pairs $(x,n)$ with $n\in\N$ and $x\in L_n$. We have
an obvious $V$-linear surjection
$$
\psi:L\to\fm
\qquad
[x,n]\mapsto xa_n.
$$
Notice that if $[x,n]\in\Ker\,\psi$, then $a_n\cdot[x,n]=[a_nx,n]=0$;
on the other hand, for every $k\geq n$ there exists $y\in V$ such that
$[x,n]=[y,k]$, whence $a_k\cdot[x,n]=a_k\cdot[y,k]=0$ as well. This
shows that $\fm\cdot\Ker\,\psi=0$. Therefore,
$\fm\otimes_V\psi:\fm\otimes_VL\to\tilde\fm:=\fm\otimes_V\fm$ is
an isomorphism. But by construction, $L$ is a flat $V$-module,
hence the natural map $\fm\otimes_VL\to\fm L$ is an isomorphism;
lastly, for every $[x,n]\in L$ we have $[x,n]=c_n[x,n+1]$, so
$\fm L=L$. Summing up, $\fm\otimes_V\psi$ is naturally identified
with the isomorphism :
$$
L\isom\tilde\fm
\qquad
[x,n]\mapsto x\otimes a_n.
$$

(iii)\ \
Furthermore, the functor $\bE$ induces a functor on almost algebras
$$
\bE^a:(A_0,\fm)^a\Alg\to(\bE_0,\fm_\bE)^a\Alg
\qquad
R^a\mapsto\bE(R/pR)^a.
$$
Indeed, let $f:R\to S$ be a homomorphism of $A_0$-algebras such
that $f^a:R^a\to S^a$ is an isomorphism of $(A_0,\fm)^a$-algebras,
and set $f_0:=f\otimes_\Z\Z/p\Z$; we need to check that
$\bE(f_0)^a:\bE(R/pR)^a\to\bE(S/pS)^a$ is an isomorphism of
$(\bE_0,\fm_\bE)^a$-algebras. However, $\bE(f_0)$ is the limit
of the system of rings $(R_n~|~n\in\N)$, with $R_n:=R/pR$ for
every $n\in\N$, and with transition map $R_{n+1}\to R_n$ given
by the Frobenius endomorphism $\Phi_{R/pR}$. Now, if $g:A_0\to R$
is the structure map of the $A_0$-algebra $R$, then $R_0$ is
naturally an $\bE_0$-algebra, via the composition $g':\bE_0\to R_0$
of the projection $\bar u_{A_0/pA_0}:\bE_0\to A_0/pA_0$ and
$g\otimes_\Z\Z/p\Z:A_0/pA_0\to R_0$. Let $\Phi_{\bE_0}$ be the
Frobenius automorphism of $\bE_0$; then for every $n\in\N$,
the ring $R_n$ is an $\bE_0$-algebra with structure morphism
$g'\circ\Phi^{-n}_{\bE_0}:\bE_0\to R_n$, {\em i.e.}
$R_n=(R_0)_{(\Phi^{-n}_{\bE_0})}$, and $\Phi_{R/pR}:R_{n+1}\to R_n$ is
a homomorphism of $\bE_0$-algebras, for these $\bE_0$-algebra
structures. Likewise, $\bE(S/pS)$ is the limit of a system
of $\bE_0$-algebras $(S_n~|~n\in\N)$, and $f$ induces a system
of maps of $\bE_0$-algebras $f_n:=f_à:R_n\to S_n$, whose inverse
limit is $\bE(f_0)$. Since $f^a$ is an isomorphism, and since
$\Phi_{\bE_0}:(\bE_0,\fm_\bE)\isom(\bE_0,\fm_\bE)$ is an isomorphism
of basic setups, it is easily seen that $f^a_n:R_n^a\to S^a_n$ is
an isomorphism of $(\bE_0,\fm_\bE)^a$-algebras, for every $n\in\N$;
on the other hand, the functor $(-)^a$ commutes with limits, since
it has a left adjoint, whence the contention.
\end{remark}

\sset\subsubsection{}\label{subsec_almost-perfectoids}
Let $A$ be a Tate ring, $A_0\subset A$ a subring of definition,
$p\in\N$ a prime integer, and $A_1\subset A$ the $p$-integral
closure of $A_0$ in $A$ (see definition
\ref{def_p-integrally-closed}(ii)).
Moreover, suppose there exists $b\in A_0^{\circ\circ}\cap A^\times$
with $pA_0\subset b^pA_0$, and let $A_0^\wedge$ and $A^\wedge_1$ be
the $b$-adic completions of $A_0$ and $A_1$. Notice that the
Frobenius endomorphism $\Phi_{A_i}$ of $A_i/b^pA_i$ induces a
ring homomorphism
$$
\bar\Phi_{A_i}:A_i/bA_i\to A_i/b^pA_i
\qquad
\text{for $i=0,1$}.
$$
Furthermore, let $(A_0,\fm)$ be a basic setup, such that $\fm$
is the filtered union of a countable system of principal subideals,
and denote by $\bar\fm$ the image of $\fm$ in $\bar A_0:=A_0/pA_0$.
Then $\fm$ is generated by the system $(x^p~|~x\in\fm)$
(\cite[Prop.2.17]{Ga-Ra}), and therefore the Frobenius endomorphism
$\Phi$ of $\bar A_0$ is a morphism of basic setups in the sense of
\cite[\S3.5]{Ga-Ra}
$$
\Phi:(\bar A_0,\bar\fm)\to(\bar A_0,\bar\fm).
$$
As detailed in \cite[\S3.5.7]{Ga-Ra}, the latter induces
a pull-back functor
$$
(\bar A_0,\bar\fm)^a\Alg\to(\bar A_0,\bar\fm)^a\Alg
\qquad
B\mapsto B_{(\Phi)}.
$$

\begin{definition}\label{def_almost-perfectoid}
In the situation of \eqref{subsec_almost-perfectoids}, we say
that the basic setup $(A_0,\fm)$ is {\em almost perfectoid} if
$A_0$ is complete and separated, and the Frobenius endomorphism
of $A_0/b^pA_0$ induces an isomorphism of $(A_0,\fm)^a$-algebras
$\bar\Phi{}^a_{A_0}:(A_0/bA_0)^a\isom(A_0/b^pA_0)^a_{(\Phi)}$.
\end{definition}

\begin{remark}
The following proposition \ref{prop_almost-perfectoid}(i) implies
in particular that definition \ref{def_almost-perfectoid} is
independent of the choice of the element $b\in A_0$ such that
$p\in b^pA_0$ : if $b'\in A_0$ is another such element, and if
the $b$-adic topology on $A_0$ agrees with the $b'$-adic topology,
then $\bar\Phi{}^a_{A_0}$ is an isomorphism if and only if the same
holds for the corresponding morphism of $(A_0,\fm)^a$-algebras
$(A_0/b'A_0)^a\isom(A_0/b'^pA_0)^a_{(\Phi)}$. See also lemma
\ref{lem_almost-perfectoid}.
\end{remark}

\begin{proposition}\label{prop_almost-perfectoid}
{\em(i)}\ \
In the situation of \eqref{subsec_almost-perfectoids}, the
following conditions are equivalent :
\begin{enumerate}
\alphaenu
\item
The basic setup $(A^\wedge_0,\fm A^\wedge_0)$ is almost perfectoid.
\item
The basic setup $(A^\wedge_1,\fm A^\wedge_1)$ is almost perfectoid,
and the inclusion map $j:A_0\to A_1$ induces isomorphisms of
$(A_0,\fm)^a$-algebras $j^a:A_0^a\isom A_1^a$.
\item
The topological closure $\bar C\subset A_0$ of the subring
$C:=\Z[a^p~|~a\in A_0]$ contains $\fm$, and $j$ induces an
isomorphism of $(A_0,\fm)^a$-algebras $j^a:A_0^a\isom A_1^a$.
\end{enumerate}

{\em(ii)}\ \
If $A_0$ is integrally closed (resp. $p$-integrally closed)
in $A$, then $A^a_{0*}$ is integrally closed (resp. $p$-integrally
closed) in $A^a_{0*}[b^{-1}]$.
\end{proposition}
\begin{proof} By simple inspection we see that
(i.b)$\Rightarrow$(i.a). Conversely, suppose that (i.a) holds,
and let also $j^\wedge:A_0^\wedge\to A_1^\wedge$ be the $b$-adic
completion of $j$; in order to show (i.b), it suffices to check :

\begin{claim}\label{cl_roads-Portishead}
$\Ker\,\bar\Phi{}^a_{A_0}=0$ $\Leftrightarrow$ $j^a$ is an
isomorphism $\Leftrightarrow$ $(j^\wedge)^a$ is an isomorphism.
\end{claim}
\begin{pfclaim} Suppose first that $\Ker\,\bar\Phi{}^a_{A_0}=0$.
We can construct $A_1$ as the increasing union of the ascending
chain $(A_{0,n}~|~n\in\N)$ of subrings of $A$ defined inductively
as follows: $A_{0,0}:=A_0$, and $A_{0,n}:=A_{0,n-1}[\Sigma_n]$, with
$\Sigma_n:=\{x\in A~|~x^p\in A_{0,n-1}\}$, for every $n>0$. Thus,
we are reduced to showing that the inclusion $A_0\to A_{0,n}$
induces an isomorphism $A_0^a\isom A^a_{0,n}$ for every $n\in\N$.
However, if for some $n>0$ the inclusion $j_{n-1}:A_{0,n-1}\to A_{0,n}$
induces an isomorphism $j^a_{n-1}:A^a_{0,n-1}\isom A^a_{0,n}$, then
the same holds for $j_n:A_{0,n}\to A_{0,n+1}$ : indeed, let
$x\in\Sigma_{n+1}$; since $x^p\in A_{0,n}$, by assumption we have
$(ax)^p\in A_{0,n-1}$ for every $a\in\fm$, whence
$ax\in\Sigma_n\subset A_{0,n}$. Thus, we are further reduced to
checking that $j_0$ induces an isomorphism $j^a_0:A_0^a\isom A_{0,1}^a$,
{\em i.e.} that $(\Coker\,j_0)^a=0$. We argue as in the proof of
lemma \ref{lem_crit-p-integr-closed}(i) : let $x\in A_1\setminus A_0$,
and denote by $m\in\N$ the smallest integer such that
$b^mx\cdot\fm\not\subset A_0$. Hence, $y:=b^{mp}x^p\in A_0$, and
$(b^{m+1}x)^p=b^py^p\in b^pA_0$. By assumption, for every $a\in\fm$
we then have $z\in A_0$ such that $ab^{m+1}x=bz$, whence
$ab^mx\in A_0$, contradicting the choice of $m$.

Conversely, suppose that $j^a$ is an isomorphism, and let
$x\in A_0$ be an element whose class $\bar x\in A_0/bA_0$
lies in $\Ker\,\bar\Phi_{A_0}$; then $x^p\in b^pA_0$, so
$(x/b)^p\in A_0$, and our assumption implies that
$\fm\cdot x\in bA_0$, {\em i.e.} $\fm\cdot\bar x=0$.

Next, if $j^a$ is an isomorphism, the same holds for
$j_k^a:=A^a_0/b^kA_0\otimes_{A_0^a}j^a:A_0^a/b^kA_0^a\to A_1^a/b^kA^a_1$,
for every $k\in\N$, and then also for the limit of the inverse
system $(j_k^a~|~k\in\N)$. But the functor $(-)^a$
commutes with limits, since it has a left adjoint, so $(j^\wedge)^a$
is an isomorphism.

Conversely, if $(j^\wedge)^a$ is an isomorphism, the induced
map $A_0^\wedge[b^{-1}]/A_0^\wedge\to A_1^\wedge[b^{-1}]/A_1^\wedge$
is an almost isomorphism of $A_0$-modules. But proposition
\ref{prop_complete-f-adic}(iii) easily implies that the
natural map $A_0[b^{-1}]/A_0\to A_0^\wedge[b^{-1}]/A_0^\wedge$
is bijective, and similarly for $A_1$; hence $j$ induces
isomorphisms $A^a_0[b^{-1}]/A^a_0\isom A^a_1[b^{-1}]/A^a_1$ and
$A_0[b^{-1}]\isom A_1[b^{-1}]$. By the 5-lemma, it follows
easily that $j^a$ is an isomorphism.
\end{pfclaim}

Next, let $\fm_{\bar C}\subset\bar C$ be the ideal generated
by $(a^p~|~a\in\fm)$. Since $\fm$ is the filtered union of
a countable system of principal subideals, the same holds
for $\bar\fm_C$; by the same token, since $\fm=\fm^2$, for
every $a\in\fm$ there exist $b,c\in\fm$ such that $a=bc$,
whence $a^p=b^pc^p$, which shows that $\fm_{\bar C}=\fm_{\bar C}^2$.
Hence, $(\bar C,\fm_{\bar C})$ is a basic setup; moreover, we
have already remarked that $\fm=\fm_{\bar C}A_0$, hence the
inclusion $\bar C\to A_0$ is a morphism of basic setups
$(\bar C,\fm_{\bar C})\to(A_0,\fm)$.

Notice then that the $\bar C$-module $A_0/\bar C$ is the limit
of the inverse system
$$
(M_n:=A_0/(b^nA_0+C)~|~n\in\N).
$$
Hence, $\fm\subset\bar C$ if and only if the
$(\bar C,\fm_{\bar C})^a$-module $(A_0/\bar C)^a$ vanishes, if
and only if $M_n^a=0$ for every $n\in\N$. Now, for given $n\in\N$,
consider the $b^p$-adic filtration on $M_n$, and let
$\gr_\bullet M_n$ be the associated graded $\bar C$-module; clearly
$\gr_iM_n$ is a quotient of $M_p$, for every $i\in\N$. Hence,
$\fm\subset\bar C$ if and only if $M_p^a=0$; the latter in turns
holds if and only if $\Coker\,\bar\Phi{}_{A_0}^a=0$. This already
shows that (a),(b)$\Rightarrow$(c). Conversely, if (c) holds,
in order to deduce (a), it suffices to check that
$\Ker\,\bar\Phi{}_{A_0}^a=0$; since $j^a$ is an isomorphism,
the latter holds by claim \ref{cl_roads-Portishead}.

(ii): Notice first that since $b\cdot\one_{A_0}$ is an
injective map, the same holds for
$(b\cdot\one_{A_0})^a_*=b\cdot\one_{A^a_{0*}}$, since the functor
$(-)^a_*$ commutes with all limits. Thus, $A^a_{0*}$ is
naturally a subring of $A^a_{0*}[b^{-1}]$. Now, let
$x\in A^a_{0*}[b^{-1}]$ be integral over $A^a_{0*}$; hence
there exist $n,k\in\N$ and $y,a_1,\dots,a_k\in A^a_{0*}$
such that $x=y/b^n$ and
$$
(y/b^n)^k+a_1\cdot(y/b^n)^{k-1}+\cdots+a_k=0.
$$
Recall that $A^a_{0*}=\Hom_{A_0}(\tilde\fm,A_0)$, with
$\tilde\fm:=\fm\otimes_{A_0}\fm$. Then we get
$(y(t)/b^n)^k+a_1(t)\cdot(y(t)/b^n)^{k-1}+\cdots+a_k(t)=0$
for every $t\in\tilde\fm$, and if $A_0$ is integrally
closed in $A$, it follows that $y(t)/b^n\in A_0$ for
every such $t$, {\em i.e.} $x\in A^a_{0*}$. One shows
likewise that if $A_0$ is $p$-integrally closed in $A$,
then $A^a_{0*}$ is $p$-integrally closed in $A^a_{0*}[b^{-1}]$.
\end{proof}

\sset\subsubsection{}\label{subsec_perf-Abhyankar}
Let $A$ be a perfectoid Tate ring, $A_0\subset A$ a perfectoid
subring of definition, and $b\in A_0^{\circ\circ}\cap A^\times$,
so that the topology of $A_0$ is $b$-adic. Let also
$g_\bullet:=(g_n~|~n\in\N)\in\bE_0:=\bE(A_0)$, and set
$g:=g_0=\bar u_{A_0}(g_\bullet)$. Hence, $g^\gamma$ is well
defined in $A_0$ for every $\gamma\in\N[1/p]$, and set
$$
\fm:=\bigcup_{n\in\N}g_nA_0
\qquad
A'_0:=\Img(A_0\to A_0[1/g])
\qquad\text{and}\qquad
A_1:=\{a\in A_0[1/g]~|~\fm\cdot a\subset A'_0\}.
$$
Clearly $(A_0,\fm)$ is a basic setup fulfilling the conditions
of \eqref{subsec_almost-perfectoids}.

\begin{lemma}\label{lem_identify-A^a_0*}
With the notation of \eqref{subsec_perf-Abhyankar}, let
$A_0^a$ be the $(A_0,\fm)^a$-algebra associated with $A_0$.
Then the natural map $A_0\to A_1$ factors through the unit
of adjunction $A_0\to A^a_{0*}$ and an isomorphism of
$A_0$-algebras $A^a_{0*}\isom A_1$.
\end{lemma}
\begin{proof} We claim that the induced surjection
$\pi:A_0\to A'_0$ induces an isomorphism of $(A_0,\fm)^a$-algebras
$\pi^a:A_0^a\isom A'^a_0$. Indeed, if $a\in\Ker\,\pi$, there exists
$n\in\N$ such that $g^na=0$ in $A_0$; then $(g^{n/p^k}a)^{p^k}=0$
in $A_0$, and consequently $g^{n/p^k}a=0$ in $A_0$ for every
$k\in\N$, due to corollary \ref{cor_perf-are-reduced}(i),
whence the contention. Since the natural map $A_0\to A_1$
factors through $\pi$ and the inclusion $i:A'_0\to A_1$, and
since $\pi^a$ is an isomorphism we are reduced to checking that
$i$ factors through the unit of adjunction $A'_0\to A'^a_{0*}$ and
an isomorphism $A'^a_{0*}\isom A_1$.

Next, we claim that the surjection $\fm\to\fm A'_0$ is bijective.
Indeed, say that $x\in\fm$ and $\pi(x)=0$; we may write
$x=g^\gamma y$ for some $\gamma\in\N[1/p]\setminus\{0\}$,
and then clearly $\pi(y)=0$. But by the foregoing, it follows
that $g^\gamma y=0$, whence the contention.

Set $\tilde\fm:=\fm\otimes_{A_0}\fm$; by the foregoing, the
isomorphism $\pi^a$ induces identifications :
$$
A^a_{0*}=\Hom_{A_0}(\tilde\fm,A_0)\isom
A'^a_{0*}=\Hom_{A_0}(\tilde\fm,A'_0)=
\Hom_{A'_0}(\fm A'_0\otimes_{A'_0}\fm A'_0,A'_0).
$$
Notice that the multiplication map
$\fm A'_0\otimes_{A'_0}\fm A'_0\to\fm A'_0$ is the colimit
of the system of isomorphisms
$(g^\gamma A'_0\otimes_{A'_0}g^\gamma A'_0\isom g^{2\gamma}A'_0~|~
\gamma\in\N[1/p]\setminus\{0\})$, and it is therefore an
isomorphism. So, finally, we get a natural identification
$$
A^a_{0*}\isom\Hom_{A'_0}(\fm A'_0,A'_0)=A_1
$$
from which the lemma follows easily : details left to the reader.
\end{proof}

The following result generalizes the ``perfectoid Abhyankar's lemma''
of \cite{YAn-2}.

\begin{theorem}\label{th_perf-Abhyankar}
In the situation of \eqref{subsec_perf-Abhyankar}, let $B$
be a finite \'etale $A_0[1/g]$-algebra, $B_1$ the integral
closure of $A_1$ in $B$, and endow $A_1$ and $B_1$ with their
$b$-adic topologies. For the almost structure given by
$(A_0,\fm)$, we have :
\begin{enumerate}
\item
The basic setup $(B_1,\!\fm B_1)$ is almost perfectoid, and $B_1$
is $p$-integrally closed in $B_1[1/b]$.
\item
$(B_1/b^nB_1)^a$ is an \'etale $(A_0/b^nA_0)^a$-algebra of finite
rank, for every $n\in\N$.
\item
If $B$ is a faithfully flat $A_0[1/g]$-algebra, then $(B_1/b^nB_1)^a$
is a faithfully flat $(A_0/b^nA_0)^a$-algebra, for every $n\in\N$.
\item
$A_1$ is complete and separated, and is integrally closed in
$A_0[1/g]$.
\item
The unit of adjunction $B_1\to B^a_{1*}$ is an isomorphism.
\end{enumerate}
\end{theorem}
\begin{proof} Arguing as in \eqref{subsec_Samuel-is-back}, we may
assume that $pA_0\subset b^pA_0$, and $b=\bar u_{A_0}(b_\bullet)$ for
some $b_\bullet\in\bE(A_0)$, so that $b^\gamma$ is well defined in $A$
for every $\gamma\in\Z[1/p]$. We let $I:=A_0b+A_0g$, and
$U:=\Spec\,A_0\setminus\Spec\,A/I$; consider the quasi-affinoid
ring $\underline A_0:=(A_0,A_0,U)$, and the rational subsets of
$X:=\Spa\,\underline A_0$
$$
R_n:=R_{A_0}\bigl(\textstyle{\frac{b^n}{g}}\bigr)\cap X
\qquad
R'_n:=R_{A_0}\bigl(\textstyle{\frac{g}{b^n}}\bigr)\cap X
\qquad
\text{for every $n\in\N$}
$$
(definition \ref{def_adic-spectrum}(iv)). Set 
$\sA_n:=\cO^{\wedge+}_X(R_n)$ and $\sA'_n:=\cO^{\wedge+}_X(R'_n)$
for every $n\in\N$; since $R_n\subset R_{n+1}$ for every $n\in\N$,
we deduce a well defined inverse system of $A_0$-algebras
$(\sA_n~|~n\in\N)$.

\begin{claim}\label{cl_ginocchio-quasi-ok}
(i)\ \
The natural map $A_0\to\cO^{\wedge+}_X(X)$ induces an isomorphism
$A_0^a\isom\cO^{\wedge+}_X(X)^a$ of $(A_0,\fm)^a$-algebras.

(ii)\ \
For every $i\in\N$, the following holds :
\begin{enumerate}
\alphaenu
\item
The induced inverse system $(\sA_n/b^i\sA_n~|~n\in\N)$ is
almost essentially constant, relative to the basic setup $(A_0,\fm)$
(see definition \ref{def_almost-Mittag-Leff}(iii)).
\item
The natural cone $(\rho_n:A_0\to\sA_n~|~n\in\N)$ induces a
universal cone $((A_0/b^iA_0)^a\to(\sA_n/b^i\sA_n)^a~|~n\in\N)$
in the category of $(A_0,\fm)^a$-algebras.
\end{enumerate}
\end{claim}
\begin{pfclaim}(i): By construction, $\cO^{\wedge+}_X(X)$ is the
integral closure of the image of $A_0$ in $A_U:=\cO_U(U)$. On
the other hand, $A^a_{0*}$ is integrally closed in $A^\circ_U$
(claim \ref{cl_perf-and-int-closed}(ii)), and recall that
$A^\circ_U$ is integrally closed in $A_U$ (remark
\ref{rem_someth-on-bdd-in-Z-lin}(iv)). It follows that the
unit of adjunction $\eta:A_0\to A^a_{0*}$ factors through the
inclusion map $\cO^{\wedge+}_X(X)\to A^a_{0*}$. Since $\eta^a$
is an isomorphism, the assertion follows.

(ii): In light of (i) and theorem \ref{th_integral-variant},
the following sequence of $(A_0,\fm)^a$-modules is short exact
for every $n\in\N$ :
$$
0\to A_0^a\to\sA^a_n\oplus\sA'^a_n\to\cO^+_X(R_n\cap R'_n)^a\to 0.
$$
Moreover, clearly the image of $b$ is a regular element in
$\cO^+_X(R_n\cap R'_n)$; by the snake lemma, we easily deduce,
for every $n,i\in\N$, a short exact sequence of $(A_0,\fm)^a$-modules
$$
0\to A^a_0/b^iA^a_0\xrightarrow{r_n}
A_0/b^iA_0\otimes_{A_0}(\sA_n\oplus\sA'_n)^a\xrightarrow{s_n}
A_0/b^iA_0\otimes_{A_0}\cO^+_X(R_n\cap R'_n)^a\to 0.
$$
Let us show that for every $\gamma\in\N[1/p]\setminus\{0\}$, every
$n\in\N$ with $n\gamma\geq i$ and every $a\in A_0$ with
$\rho_n(a)\in b^i\sA_n$, we have $g^{2\gamma}a\in b^iA_0$. Indeed,
let $\bar a\in A_0/b^iA_0$ be the image of $a$; then
$r_n(\bar a)=(0,\bar x)$ for some $\bar x\in\sA'_n$, and we
notice that $g^\gamma\cdot\bar x=b^{n\gamma}\cdot(g/b^n)^\gamma\cdot\bar x$,
whence $g^\gamma\cdot\bar x=0$, since $(g/b^n)^\gamma\in\sA'_n$.
It follows that $r_n(g^\gamma\cdot\bar a)=0$, whence
$g^{2\gamma}\cdot\bar a=0$, as required. Lastly, let us check
that for every $\gamma\in\N[1/p]\setminus\{0\}$, every
$n\in\N$ with $n\gamma\geq i$ and every $\bar x\in\sA_n/b^i\sA_n$
there exists $\bar a\in A_0/b^iA_0$ with
$\rho_n(\bar a)=g^{2\gamma}\cdot\bar x$. To this aim, set
$\bar y:=s_n(\bar x,0)$; again, we notice that
$g^\gamma\cdot\bar y=b^{n\gamma}\cdot(g/b^n)^\gamma\cdot\bar y=0$,
so that $(g^\gamma\bar x,0)\in\Ker\,s_n=\Img\,r_n$, whence the
contention. This completes the proof of (ii.a). Assertion (ii.b)
follows from the proof of (ii.a) and lemma
\ref{lem_a.ess.0-is-Serre}(iii).
\end{pfclaim}

Recall that $\cO^\wedge_X(R_n)=\sA_n[1/g]$, and set
$\sB_n:=\sA_n[1/g]\otimes_{A_0}B$ for every $n\in\N$; recall
as well that $\sA_n$ is perfectoid for its $b$-adic topology
(proposition \ref{prop_Letta-non-caduto}(iii) and theorem
\ref{th_int-subrings-perfectoid}(iii)).
Since the functor $(-)^a$ commutes with limits, from claim
\ref{cl_ginocchio-quasi-ok} we deduce that also the cone
$(\rho_n~|~n\in\N)$ is universal, {\em i.e.} we have an
isomorphism of $(A_0,\fm)^a$-algebras :
\set\begin{equation}\label{eq_there-is-no-+}
A_0^a\isom\lim_{i\in\N}\lim_{n\in\N}\sA^a_n/b^i\sA^a_n
\isom\lim_{n\in\N}\sA^a_n.
\end{equation}

\begin{claim}\label{cl_speriamo}
$\fm\sA_n$ is an open radical ideal of $\sA_n$, for every $n\in\N$.
\end{claim}
\begin{pfclaim} Clearly $\fm\sA_n$ is an open ideal of
$\sA_n$, since $(b^n/g)^\gamma\in\sA_n$ for every
$\gamma\in\N[1/p]$. Moreover, set
$J:=\bigcup_{\gamma\in\N[1/p]\setminus\{0\}}b^\gamma\sA_n$;
by corollary \ref{cor_BDS}, the quotient $\sA_n/J$ is
a perfect $\F_p$-algebra; it follows easily that $\fm/J$
is a radical ideal of $\sA_n/J$, whence the claim.
\end{pfclaim}

Say that $B$ is a projective $A_0[1/g]$-module of rank $\leq r$;
from claim \ref{cl_speriamo} and theorem
\ref{th-alm-purity-form-perfectoid}, we deduce that for every
$n\in\N$ there exists an \'etale $\sA_n^a$-algebra $C_n$ of
rank $\leq r$, unique up to unique isomorphism, with an isomorphism
of $\sA_n[1/g]$-algebras $C_n[1/g]\isom\sB_n$. From the uniqueness
property of $C_n$, there follows a unique
isomorphism of $\sA_n^a$-algebras :
\set\begin{equation}\label{eq_denote-by-D}
\sA^a_n\otimes_{\sA^a_{n+1}}C_{n+1}\isom C_n
\qquad
\text{for every $n\in\N$}.
\end{equation}
Denote by $D$ (resp. by $E$) the inverse limit of the resulting
system $(C_n~|~n\in\N)$ of $A_0^a$-algebras (resp. $(\sB_n~|~n\in\N)$
of $A_0$-algebras), and for every $i\in\N$ let as well $D_i$ be
the limit of the induced system $(C_n/b^iC_n~|~n\in\N)$ of
$A^a_0/b^iA^a_0$-algebras. By proposition
\ref{prop_return-of-AlMod}(ii,iii), $(C_n/b^iC_n~|~n\in\N)$ is
almost essentially constant, and the projection $D_i\to C_n/b^iC_n$
induces an isomorphism of $A^a_0/b^iA^a_0$-algebras :
\set\begin{equation}\label{eq_never-too-late}
\sA^a_n\otimes_{A^a_0}D_i\isom C_n/b^iC_n
\qquad
\text{for every $n\in\N$}.
\end{equation}
Notice that $b$ and $g$ are regular elements of $\sA_n$, and
then they are also regular elements of $C_{n*}$, since $C_n$ is
a flat $\sA^a_n$-algebra; in view of the induced isomorphism
\set\begin{equation}\label{eq_call-of-hagar}
D_*\isom\lim_{n\in\N}C_{n*}
\end{equation}
(proposition \ref{prop_was-get-maddd}(iii)) it follows easily
that $b$ and $g$ are regular elements of $D_*$.

\begin{claim}\label{cl_Trump-in-UK}
(i)\ \
For every $n\in\N$, the ring $C_{n*}$ is perfectoid
for its $b$-adic topology.

(ii)\ \
$D_*$ is $p$-integrally closed in $D_*[1/b]$ and is integrally
closed in $D_*[1/g]$.

(iii)\ \
The $A^a_0/b^iA_0$-algebra $D_i$ is \'etale of rank $\leq r$
for every $i\in\N$.

(iv)\ \
If $B$ is a faithfully flat $A_0[1/g]$-algebra, then $D_i$
is a faithfully flat $A^a_0/b^iA_0$-algebra for every $i\in\N$.

(v)\ \
$B$ and $D_*$ are integrally closed in $E$.

(vi)\ \
The basic setup $(D_*,\fm D_*)$ is almost perfectoid for the
$b$-adic topology of $D_*$, and the projection $D\to D_i$
factors through an isomorphism $D/b^iD\isom D_i$ for every $i\in\N$.
\end{claim}
\begin{pfclaim} (i) follows from proposition
\ref{prop_alm-etale-over-perf-is-perf} and claim \ref{cl_speriamo},
since we have already recalled that $\sA_n$ is perfectoid for its
$b$-adic topology. Together with corollary \ref{cor_variant-BDS}
and lemma \ref{lem_crit-p-integr-closed}(i) we deduce that $C_{n*}$
is $p$-integrally closed in $C_{n*}[1/b]$, for every $n\in\N$.
Then lemma \ref{lem_crit-p-integr-closed}(iii.a) and the isomorphism
\eqref{eq_call-of-hagar} imply that $D_*$ is $p$-integrally closed
in the limit $L$ of the system $(C_{n,*}[1/b]~|~n\in\N)$. However,
the induced map $D_*\to L$ factors through an injective map
$D_*[1/b]\to L$, so $D_*$ is also $p$-integrally closed in $D_*[1/b]$.

Similarly, $C_{n*}$ is integrally closed in $C_{n*}[1/g]=B_n$ for
every $n\in\N$, due to claim \ref{cl_perf-and-int-closed}(ii);
it follows easily that $D_*$ is integrally closed in $E$. But the
induced map $D_*\to E$ factors through an injective map
$D_*[1/g]\to E$, so $D_*$ is also integrally closed in $D_*[1/g]$.

(iii) follows from theorem \ref{th_alm-ess-cnst-inv-sys}(iii,iv).

(iv): If $B$ is a faithfully flat $A_0[1/g]$-algebra, $C_n$ is
a faithfully flat $\sA_n^a$-algebra for every $n\in\N$ (theorem
\ref{th-alm-purity-form-perfectoid}(ii)); then the assertion
follows from theorem \ref{th_alm-ess-cnst-inv-sys}(ii)).

(v): We have already observed that $D_*$ is integrally closed
in $E$. Next, by construction, $\sB_n=B\otimes_{A_0[1/g]}\sA_n[1/g]$
for every $n\in\N$. Since $B$ is a projective $A_0[1/g]$-module
of finite rank, there follows a natural identification :
$$
E\isom B\otimes_{A_0[1/g]}\lim_{n\in\N}\sA_n[1/g].
$$
Since $B$ is an \'etale $A_0[1/g]$-algebra, by proposition
\ref{prop_8231} we are then reduced to showing that $A_0[1/g]$
is integrally closed in the limit $L$ of the system of
$A_0$-algebras $(\sA_n[1/g]~|~n\in\N)$. However, invoking
again claim \ref{cl_perf-and-int-closed}(ii) we see that
$\sA^a_{n*}$ is integrally closed in $\sA_n[1/g]$ for every
$n\in\N$, hence the limit $L'$ of the system $(\sA^a_{n*}~|~n\in\N)$
is integrally closed in $L$, and therefore $L'[1/g]$ is integrally
closed in $L$. But in view of \eqref{eq_there-is-no-+}, we have
a natural identification of $A_0$-algebras :
$A_0[1/g]\isom A^a_{0*}[1/g]\isom L'[1/g]$, whence the assertion.

(vi):  Since $C_n$ is the limit of the inverse system
$(C_n/b^iC_n~|~i\in\N)$ for every $n\in\N$, we obtain a natural
isomorphism of $A^a_0$-algebras :
$$
D\isom\lim_{i\in\N}D_i.
$$
On the other hand, in view of \eqref{eq_never-too-late} we have
a natural isomorphism of $A^a_0$-algebras
$$
(D_{i+1}\otimes_{A^a_0}A_0^a/b^iA^a_0)\otimes_{A^a_0}\sA^a_n\isom
C_n/b^iC_n
\qquad
\text{for every $n,i\in\N$}
$$
whence, again by theorem \ref{th_alm-ess-cnst-inv-sys}(i), an
induced isomorphism of $A^a_0$-algebras :
\set\begin{equation}\label{eq_full-details}
D_{i+1}\otimes_{A^a_0}A_0^a/b^iA^a_0\isom D_i
\qquad
\text{for every $i\in\N$}.
\end{equation}
For every $i\in\N$, let $\pi_i:D\to D_i$ be the projection;
notice that $I:=\Ker\,\pi_1=bD+\Ker\,\pi_{i+1}$ for every $i\in\N$,
in light of the isomorphism $D_{i+1}/bD_{i+1}\isom D_1$ deduced
from \eqref{eq_full-details}. By \cite[Lemma 5.3.5(i)]{Ga-Ra},
it follows that $I=\bigcap_{i\in\N}(bD+\bar{I^{i+1}})=\bar{bD}$
(where $\bar J$ denotes the closure of $J$, for every ideal
$J\subset D$ : see \cite[Def.5.3.1(iii)]{Ga-Ra}). By
\cite[Lemma 5.3.8(iii)]{Ga-Ra} we deduce that $I=bD$, and more
generally $b^iD=\Ker\,\pi_i$ for every $i\in\N$. On the other
hand, $\pi_i$ is an epimorphism for every $i\in\N$, by lemma
\ref{lem_lim-one-special-case}(i), so it induces an isomorphism
$D/b^iD\isom D_i$, for every $i\in\N$. Hence, the natural morphism
$D\to\lim_{i\in\N}D/b^iD$ is an isomorphism of $A^a_0$-algebras, so
we have an induced isomorphism of $A_1$-algebras
$D_*\isom\lim_{i\in\N}(D/b^iD)_*$ (lemma \ref{lem_identify-A^a_0*}).
As we have already noticed, $b$ is regular in $D_*$, so for every
$i\in\N$ we get a short exact sequence of $A^a_0$-modules
$0\to D\xrightarrow{b^i\cdot\one_{D_*}}D\to D/b^iD\to 0$
and since the functor $(-)_*$ is left exact, we deduce
an injective $A_1$-linear map
$$
D_*/b^iD_*\to(D/b^iD)_*
\qquad
\text{for every $i\in\N$}
$$
from which it follows easily that also the induced map
$D_*\to\lim_{i\in\N}D_*/b^iD_*$ is an isomorphism, hence $D_*$
is $b$-adically complete and separated. In order to show that
$(D_*,\fm D_*)$ is almost perfectoid, it remains only to check
that the Frobenius endomorphism of $D_*/pD_*$ induces an isomorphism
$D/bD\isom(D/b^pD)_{(\Phi)}$ of $A^a_0/bA^a_0$-algebras (where
$\Phi$ denotes the Frobenius endomorphism of $A_0/pA_0$).
But we have already seen that $D/b^iD$ is the inverse limit
of the system $(C_n/b^iC_n~|~n\in\N)$, for every $i\in\N$,
hence we are reduced to checking that the Frobenius endomorphism
of $C_{n*}/pC_{n*}$ induces an isomorphism
$C_n/bC_n\isom(C_n/b^pC_n)_{(\Phi)}$ of $A^a_0/bA^a_0$-algebras,
for every $n\in\N$. The latter is clear, since we know already
that $C_{n*}$ is perfectoid (corollary \ref{cor_jackob}).
\end{pfclaim}

Next, recall that for every $n\in\N$ we have $A^a_0$-algebras
$\sA_{n,0},\dots,\sA_{n,r}$, and an isomorphism of $A^a_0$-algebras
$\sA_n^a\isom\sA_{n,0}\times\cdots\times\sA_{n,r}$ such that
$C_{n,k}:=\sA_{n,k}\otimes_{\sA^a_n}C_n$ is an $\sA_{n,k}$-module
of constant rank $k$, for every $k=0,\dots,r$
(\cite[Prop.4.3.27]{Ga-Ra}). For every $n\in\N$ and
$k=0,\dots,r$ consider the map of sets :
$$
\chi^{(n,k)}:C_{n,k*}\to\sA^a_{n,k*}[X]
\qquad
c\mapsto\chi_c(X)
$$
which associates with every $c\in C_{n,k*}$ the characteristic
polynomial $\chi_c(X)$ of the $\sA_{n,k}^a$-linear endomorphism
$c\cdot\one_{C_{n,k}}$ of $C_{n,k}$ (see \cite[\S4.4.29]{Ga-Ra}); by
\cite[Prop.4.4.30]{Ga-Ra} we have $\chi_c(c)=0$ for every such
$c$. By virtue of \eqref{eq_denote-by-D}, we see that the restriction
map $\sA_{n+1}\to\sA_n$ induces isomorphisms of $A^a_0$-algebras
$$
\sA^a_n\otimes_{\sA^a_{n+1}}\sA_{n+1,k}\isom\sA_{n,k}
\qquad
\text{for every $k=0,\dots,r$}
$$
and from the compatibility of traces with ring extensions
(\cite[Prop.4.1.8(ii)]{Ga-Ra}), we deduce for every
$n\in\N$ and $k=0,\dots,r$ a commutative diagram :
$$
\xymatrix@C+20pt{
C_{n+1,k*} \ar[r]^-{\chi^{(n+1,k)}} \ar[d] & \sA^a_{n+1,k*}[X] \ar[d] \\
C_{n,k*} \ar[r]^-{\chi^{(n,k)}} & \sA^a_{n,k*}[X]
}$$
whose vertical arrows are induced by the induced projection
$\sA_{n+1,k}\to\sA_{n,k}$. For every $k=0,\dots,r$, let $\sA^{(k)}$
be the limit of the induced inverse system of $A^a_0$-algebras
$(\sA_{n,k}~|~n\in\N)$, and set $D^{(k)}:=\sA^{(k)}\otimes_{A_0^a}D$;
in view of \eqref{eq_there-is-no-+}, there follows an isomorphism
of $A^a_0$-algebras $D\isom D^{(0)}\times\cdots\times D^{(r)}$.
Notice as well that $\chi_c(X)$ is a monic polynomial of degree
$k$, for every $n\in\N$, every $k=0,\dots,r$ and and every
$c\in C_{n,k*}$; in light of \eqref{eq_call-of-hagar}, it follows
easily that the limit of the system of maps $(\chi^{(n,k)}~|~n\in\N)$
is a well defined map of sets
$$
\chi^{(k)}:D^{(k)}_*\to\sA^{(k)}_*[X]
\qquad
d\mapsto\chi_d(X)
\qquad
\text{for every $k=0,\dots,r$}.
$$
Moreover, $\chi_d(X)$ is a monic polynomial with $\chi_d(d)=0$
for every $k=0,\dots,r$ and every $d\in D^{(k)}_*$. Taking into
account lemma \ref{lem_identify-A^a_0*}, it follows that $D_*$
is an integral $A_1$-algebra. Since both $D_*$ and $B$ are
integrally closed subrings of $E$ (claim \ref{cl_Trump-in-UK}(v)),
and since $A_1\subset A_0[1/g]$, we deduce that $D_*\subset B$,
and hence $D_*$ is also the integral closure of $A_1$ in $B$,
{\em i.e.} $D_*=B_1$. Then assertions (i), (ii) and (iii) of the
theorem follow from claim \ref{cl_Trump-in-UK}.

(iv): By construction, $A_1$ is an open bounded subring of $A$,
hence it is complete and separated for the topology induced by
$A$, which agrees with the $b$-adic topology (proposition
\ref{prop_f-adics}(i,ii)). In order to show that $A_1$ is
integrally closed in $A_0[1/g]$, consider the case where
$B=A_0[1/g]$; then $C_n=A_{0,n}$ for every $n\in\N$, so that
$B_1=D_*=A^a_{0*}=A_1$, whence the assertion.
\end{proof}

\begin{remark}
(i)\ \
In the situation of theorem \ref{th_perf-Abhyankar}, let
$g'_\bullet:=(g'_n~|~n\in\N)\in\bE_0$ be another element,
and $g':=g'_0=\bar u_{A_0}(g'_\bullet)$. Let also
$A_0'':=\Img(A_0\to A_0[1/(gg')])$ and set
$$
\fm':=\bigcup_{n\in\N}g'_nA_0
\qquad
\fm'':=\fm\cdot\fm'
\qquad
A_2:=\{a\in A_0[1/(gg')]~|~\fm''\cdot a\subset A_0''\}.
$$
Clearly $(A_0,\fm')$ and $(A_0,\fm'')$ are two other basic setups;
for every $A_0$-module $M$, let $(M,\fm)^a\in(A_0,\fm)^a\Mod$ (resp.
$(M,\fm')^a\in(A_0,\fm')^a\Mod$, resp.
$(M,\fm'')^a\in(A_0,\fm'')^a\Mod$) be the image of $M$.
We may apply theorem \ref{th_perf-Abhyankar} with $(A_0,\fm)$
and $B$ replaced respectively by $(A_0,\fm'')$, and $B[1/g']$ :
hence, let $B_2$ be the integral closure of $A_2$ in $B[1/g']$;
then $B_2=(B_2,\fm'')^a_*$ is an almost perfectoid $A_0$-algebra
relative to $(A_0,\fm'')$, and $(B_2/b^nB_2,\fm'')^a$ is an \'etale
$(A_0/b^nA_0,\fm'')^a$-algebra of finite rank, for every $n\in\N$.
The localizations $A_0[1/g]\to A_0[1/(gg')]$ and $B\to B[1/g']$
clearly restrict to ring homomorphisms $A_1\to A_2$ and $B_1\to B_2$.

Now, $(A_0,\fm'')^a$ is integrally closed in $(A_0[1/(gg')],\fm'')^a$
(theorem \ref{th_perf-Abhyankar}(iv) and lemma
\ref{lem_identify-A^a_0*}), hence the same holds for
$(A_0[1/g],\fm'')^a$; since $B$ is a finite \'etale
$A_0[1/g]$-algebra, it follows that $(B,\fm'')^a$ is integrally
closed in $(B[1/g'],\fm'')^a$ (proposition \ref{prop_8231}).
Taking into account theorem \ref{th_perf-Abhyankar}(v), remark
\ref{rem_prod_of-alm-structures}(ii) and \cite[Lemma 8.2.28]{Ga-Ra}
we deduce :
$$
\begin{aligned}
B_2=(B_2,\fm'')^a_*
&=\mathrm{i.c.}((A_0,\fm'')^a,(B[1/g'],\fm'')^a)_* \\
&=\mathrm{i.c.}((A_0,\fm'')^a,(B,\fm'')^a)_* \\
&=(\mathrm{i.c.}((A_0,\fm'')^a_*,(B,\fm'')^a_*),\fm'')^a_* \\
&=(\mathrm{i.c.}((A_1,\fm')^a_*,(B,\fm')^a_*),\fm'')^a_* \\
&=((\mathrm{i.c.}(A_1,B),\fm')^a_*,\fm'')^a_* \\
&=((B_1,\fm')^a_*,\fm'')^a_* \\
&=(B_1,\fm'')^a_* \\
&=(B_1,\fm')^a_*.
\end{aligned}
$$

(ii)\ \
As a special case of (i), we may take $g'_\bullet:=b_\bullet$ such
that $b:=b_0\in A^\times\cap A_0^{\circ\circ}$.
Then, endow $B_1$ with its $b$-adic topology, so that $B_1$ is
a subring of definition of the Tate ring $B_1[1/b]$, and by
theorem \ref{th_perf-Abhyankar}(i) it is $p$-integrally closed
in $B_1[1/b]$, hence $B_1[1/b]^{\circ\circ}\subset B_1$ (example
\ref{ex_one-uniform-example}). It follows easily that
$(B_1[1/b]^\circ,\fm')^a=(B_1,\fm')^a$, so we have inclusions
of rings :
$$
B_1\subset B^\circ\subset(B_1,\fm')^a_*=
(B_1[1/b]^\circ,\fm')^a_*\subset B_1[1/b].
$$
On the other hand, it is easily seen that
$(B_1,\fm')^a_*\subset B_1[1/b]^\circ$, so
$B_2=(B_1,\fm')^a_*=B_1[1/b]^\circ$.

(iii)\ \
In the situation of theorem \ref{th_perf-Abhyankar}, suppose
that $A$ is an $\F_p$-algebra; then $A_0=\bE_0$ is a perfect
$\F_p$-algebra, and it follows easily that the same holds for
$A_1$. Then, set $X:=\Spec\,A_1$ and $Z:=\Spec\,A_1/gA_1$; by
theorem \ref{th_perfect-purity}, the pair $(X,Z)$ is normal
and almost pure, hence $B^a_1$ is an \'etale $A^a_1$-algebra
of finite rank (proposition \ref{prop_almost-pure-crit}).

(iv)\ \
If $A$ is not an $\F_p$-algebra, we know neither whether $A_1$
is perfectoid, nor whether $B^a_1$ is an \'etale $A^a_1$-algebra.
However, by theorem \ref{th_perf-Abhyankar}(ii) and lemma
\ref{lem_axiomatize-familiar}(ii), for every $n\in\N$ we have
a well defined map
$$
\chi_{B_1/b^nB_1}:(B_1/b^nB_1)^a_*\to(A_0/b^nA_0)^a_*[T]
\qquad
b\mapsto\det((1+bT)\cdot\one_{B^a_1/b^nB^a_1})
$$
and clearly the inverse system of maps $(\chi_{B_1/b^nB_1}~|~n\in\N)$
amounts to a well defined map
$$
\chi_{B_1}:B_1\to A_1[T].
$$
By the same token, we may associate to the finite \'etale
$A_0[1/g]$-algebra $B$ the map
$$
\chi_B:B\to A_0[1/g,T]
\qquad
b\mapsto\det((1+bT)\cdot\one_B)
$$
and we claim that
$$
\chi_B(b)=\chi_{B_1}(b)
\qquad
\text{for every $b\in B_1$}.
$$
For the proof, define the inverse system of $A_0$-algebras
$\sA_\bullet:=(\sA_n~|~n\in\N)$, and for every $n\in\N$ the
\'etale $\sA_n^a$-algebra $C_n$ of rank $\leq r$, as in the
proof of theorem \ref{th_perf-Abhyankar}; recall that $A_1$
(resp. $B_1$) is the limit of $\sA_\bullet$ (resp. of
$(C_{n*}~|~n\in\N)$). It follows easily that $\chi_{B_1}$
is also the inverse limit of the system of similar maps
$(\chi_{C_n}~|~n\in\N)$. On the other hand, $A_0[1/g]$ is
integrally closed in the limit $L$ of the induced inverse system
$(\sA_n[1/g]~|~n\in\N)$, and $B$ is the integral closure of
$A_0[1/g]$ in the inverse system $E$ of the induced inverse system
$(C_n[1/g]~|~n\in\N)$; also, the natural map $E\to B\otimes_{A_0[1/g]}L$
is an isomorphism (see the proof of claim \ref{cl_Trump-in-UK}).
Thus, both $\chi_B$ and $\chi_{B_1}$ are restrictions of the
corresponding map $\chi_E:E\to L$, whence the assertion.

(v)\ \
Furthermore, set $R:=B_1\otimes_{A_1}B_1$, let $R^\wedge$ be the
$b$-adic completion of $R$, and
$$
e_{B/A_0[1/g]}\in R[1/g]=B\otimes_{A_0[1/g]}B
$$
the diagonal idempotent of the $A_0[1/g]$-algebra $B$; we notice :
\end{remark}

\begin{corollary}\label{cor_perf-Abhyankar}
In the situation of theorem {\em\ref{th_perf-Abhyankar}}, we have :

{\em(i)}\ \
The following conditions are equivalent:
\begin{enumerate}
\alphaenu
\item
$g^\gamma\cdot e_{B/A_0[1/g]}\in R':=\Img(R\to R[1/g])$
for every $\gamma\in\N[1/p]\setminus\{0\}$.
\item
For every $\gamma\in\N[1/p]\setminus\{0\}$ there exists
$N_\gamma\in\N$ such that
$b^{N_\gamma}g^\gamma\cdot e_{B/A_0[1/g]}\in R'$.
\item
$B^a_1$ is an \'etale $A^a_1$-algebra.
\item
$B^a_1$ is a weakly unramified $A^a_1$-algebra.
\item
$B^a_1$ is an almost finitely generated $A^a_1$-module.
\item
$B^a_1$ is an almost finite projective $A^a_1$-module.
\item
$B^a_1[1/b]$ is a weakly unramified $A^a_1[1/b]$-algebra.
\item
$B^a_1[1/b]$ is an almost finitely generated $A^a_1[1/b]$-module.
\end{enumerate}

{\em(ii)}\ \
There exists $n\in\N$ such that $g^n$ annihilates the kernel
and cokernel of the completion map $R\to R^\wedge$. Also,
$\Ann_{R^\wedge}(g)^a=0$.
\end{corollary}
\begin{proof} Obviously (c)$\Rightarrow$(a),(d); also
(d)$\Rightarrow$(g), (a)$\Rightarrow$(b) and
(f)$\Rightarrow$(e)$\Rightarrow$(h).

(a)$\Rightarrow$(c),(f): The assumption implies that $e_{B/A_0[1/g]}$
lies in the image of the induced map
$$
R^a_*\to B\otimes_{A_0[1/g]}B.
$$
Then the assertion follows from proposition
\ref{prop_axiomatize-familiar} and theorem
\ref{th_perf-Abhyankar}(iv)).

(f)$\Rightarrow$(c): We consider the commutative diagram
of $A^a_1$-algebras :
\set\begin{equation}\label{eq_very-late}
{\diagram R^a \ar[r]^-\phi \ar[d]_{\mu^a} &
\lim_{n\in\N}R^a/b^nR^a \ar[d] \\
B^a_1 \ar[r] & \lim_{n\in\N}B^a_1/b^nB^a_1
\enddiagram}
\end{equation}
where $\mu^a$ is the multiplication law, and the right vertical
arrow is the inverse limit of the system of multiplication laws
$(R^a/b^nR\to B^a_1/b^nB^a_1~|~n\in\N)$.
According to theorem \ref{th_perf-Abhyankar}(i), $B_1$ is
$b$-adically complete and separated; then by lemma
\ref{lem_lim-one-special-case}(iii) the morphism $\phi$ is an
isomorphism. By theorem \ref{th_perf-Abhyankar}(ii), for each
$n\in\N$ we have a diagonal idempotent $e_n\in(R/b^nR)_*$, and
it is easily seen that the system $(e_n~|~n\in\N)$ corresponds,
under the isomorphism $\phi$, to an idempotent element
$e_{B^a_1/A^a_1}\in R^a_*$ such that $\mu(_{B^a_1/A^a_1})=1$ and
$x\cdot e_{B^a_1/A^a_1}=0$ for every $x\in(\Ker\,\mu^a)_*$. Then
the assertion follows from \cite[Prop.3.1.4]{Ga-Ra}.

(d)$\Rightarrow$(a): Let $J:=\bigcup_{n\in\N}\Ann_R(g^n)$ and
$\bar R:=R/J$. By construction, $g$ is a regular element of
$B_1$, so that $JB_1=0$, for the $R$-algebra structure on $B_1$
given by the multiplication map $\mu:R\to B_1$. On the other
hand, by assumption $\mu^a$ is a flat morphism of $A^a_1$-algebras;
hence $\mu^a\otimes_R\bar R:\bar R\to B_1$ is still flat. Also,
the localisation $R\to R[1/g]$ factors through a monomorphism
$\bar R\to R[1/g]$ of $A^a_1$-algebras, and since $B$ is an
\'etale $A_0[1/g]$-algebra, $B=B_1[1/g]$ is a projective
$R[1/g]$-module of finite rank.
By \cite[Prop.2.4.19]{Ga-Ra} it follows that $B_1$ is an almost
finite projective $\bar R$-module, and therefore
$\Ker\,(\mu^a\otimes_R\bar R)$ is generated by an idempotent almost
element $d\in\bar R_*$ (\cite[Rem.3.1.8]{Ga-Ra}). For every
$\gamma\in\N[1/p]\setminus\{0\}$, we then find $e_\gamma\in R$
whose image in $(R/J)_*$ agrees with $g^\gamma\cdot(1-d)$. It
follows easily that the image of $e_\gamma$ in $R[1/g]$ agrees
with $g^\gamma e_{B/A_0[1/g]}$, whence the contention.

(g)$\Rightarrow$(b): The assumption means that
$\mu^a\otimes_{A_1}A_1[1/b]:R[1/b]^a\to B_1[1/b]^a$ is a flat
morphism, hence the same holds for $\mu^a\otimes_R\bar R[1/b]$,
and arguing as in the foregoing we find for every
$\gamma\in\N[1/p]\setminus\{0\}$ an integer $N_\gamma\in\N$
and $e_\gamma\in R$ whose image in $R[1/(gb)]$ agrees with the
image of $g^\gamma b^{N_\gamma}e_{B/A_0[1/g]}$. Notice that, since
$B_1[1/g]=B$ is a flat $A_1[1/g]$-algebra, and since $b$
is regular in $A_1$, it follows that $b$ is regular also
in $R[1/g]$, hence the localization $R[1/g]\to R[1/(gb)]$
is injective, and we deduce that the image of $e_\gamma$ agrees
with $g^\gamma b^{N_\gamma}e_{B/A_0[1/g]}$ {\em already in} $R[1/g]$,
as required.

(b)$\Rightarrow$(a): Define the inverse system of $A_0$-algebras
$(\sA_n~|~n\in\N)$ as in the proof of theorem
\ref{th_perf-Abhyankar}, and for every $n\in\N$ let $C_n$ be the
unique \'etale $\sA_n^a$-algebra up to unique isomorphism, such
that $C_n[1/g]$ is isomorphic to $\sA_n[1/g]\otimes_{A_0}B$.
By claim \ref{cl_ginocchio-quasi-ok}(ii), for every $i\in\N$
the induced inverse system $(\sA_n/b^i\sA_n~|~n\in\N)$ is almost
essentially constant relative to the basic setup $(A_0,\fm)$.
In view of \eqref{eq_denote-by-D} the same follows for the
inverse system $(C^a_n/b^iC^a_n~|~n\in\N)$, for every $i\in\N$;
but the proof of theorem \ref{th_perf-Abhyankar} also shows
that the inverse limit of the latter system is naturally
isomorphic to $B^a_1/b^iB^a_1$, for every $i\in\N$. Thus, set
$R_n:=C_n\otimes_{\sA_n^a}C_n$ for every $n\in\N$; we deduce that
for every $i\in\N$ and there exists $n(i,\gamma)\in\N$ such that
$g^\gamma$ annihilates the kernel of the induced morphism
$$
R^a/b^iR^a\to R_{n(i,\gamma)}/b^iR_{n(i,\gamma)}.
$$
Notice also that, since $g$ is regular in $\sA_n$, and $C_n$ is
a flat $\sA_n^a$-algebra, then $g$ is regular in $C_n$, and the
localization $R_n\to R_n[1/g]$ is a monomorphism of
$\sA^a_n$-modules, for every $n\in\N$. Moreover, the image of
$e_{B/A_0[1/g]}$ in $R_n[1/g]$ is the diagonal idempotent of the
\'etale $\sA_n[1/g]$-algebra $C_n[1/g]$, and the latter agrees
with the image of the diagonal idempotent $e_n\in R_{n*}$ of the
\'etale $\sA^a_n$-algebra $C_n$. Summing up, we conclude that for
every $\gamma\in\N[1/p]\setminus\{0\}$, the image of $e_\gamma$
vanishes in $R_{n(N_\gamma,\beta)}/b^{N_\gamma}R_{n(N_\gamma,\gamma)}$, and
hence $g^{2\gamma}e_\gamma$ vanishes in $R/b^{N_\gamma}R$, {\em i.e.}
$g^{2\gamma}e_\gamma=b^{N_\gamma}e'_\gamma$ for some $e'_\gamma\in R$.
Since $b$ is regular in $R[1/g]$, it follows that the image of
$e'_\gamma$ in $R[1/g]$ agrees with $g^{3\gamma}e_{B/A_0[1/g]}$, as required.

(h)$\Rightarrow$(b): By assumption, for every
$\gamma\in\N[1/p]\setminus\{0\}$ there exists $n(\gamma)\in\N$
and an $A_1$-linear map $f_\gamma:A^{\oplus n(\gamma)}_1\to B_1$ such
that
\set\begin{equation}\label{eq_banacchia}
g^\gamma\cdot\Coker(A_1[1/b]\otimes_{A_1}f_\gamma)=0.
\end{equation}
Then, let $E_\gamma$ be the fibre product in the cartesian diagram
of $A_1$-modules :
$$
\cD_\gamma \qquad : \qquad
{\diagram E_\gamma \ar[r]^-{f'_\gamma} \ar[d] &
B_1 \ar[d]^{g^\gamma\one_{B_1}} \\
A_1^{\oplus n(\gamma)} \ar[r]^-{f_\gamma} & B_1.
\enddiagram}$$
By virtue of \eqref{eq_banacchia}, the image of
$A_1[1/b]\otimes_{A_1}g^\gamma\one_{B_1}$ lies in the image of
$A_1[1/b]\otimes_{A_1}f_\gamma$; on the other hand, the diagram
$A_1[1/b]\otimes_{A_1}\cD_\gamma$ is still cartesian, so the map
$A_1[1/b]\otimes_{A_1}f'_\gamma$ {\em is surjective}. Recall that
$b$ is regular in both $A_1$ and $B_1$, therefore also in
$E_\gamma$; then endow $A_1[1/b]$ (resp. $B_1[1/b]$) with the f-adic
topology $\cT_{A_1}$ (resp. $\cT_{B_1}$) such that $A_1$ (resp.
$B_1$) is a subring of definition, with ideal of definition
$bA_1$ (resp. $b_1B_1$); since $A_1$ and $B_1$ are $b$-adically
complete and separated (theorem \ref{th_perf-Abhyankar}(i)),
then $\cT_{A_1}$ and $\cT_{B_1}$ are complete and separated
(proposition \ref{prop_complete-f-adic}). Endow also
$A_1[1/b]^{\oplus n(\gamma)}$ with the product topology
$\cT_{A_1}\times\cdots\times\cT_{A_1}$; it is easily seen that
all the maps of $A_1[1/b]\otimes_{A_1}\cD_\gamma$ are continuous
for these topologies, hence $A_1[1/b]\otimes_{A_1}E_\gamma$
is a closed subset of
$A_1[1/b]\otimes_{A_1}(A_1^{\oplus n(\gamma)}\oplus B_1)$, and
is therefore complete and separated for the topology induced
by the inclusion in this $A_1$-module; moreover, since
$A_1^{\oplus n(\gamma)}$ is open in $A_1[1/b]^{\oplus n(\gamma)}$
and $B_1$ is open in $B_1[1/b]$, we see that $E_\gamma$ {\em is
open in} $A_1[1/b]\otimes_{A_1}E_\gamma$. By the Banach open
mapping theorem (\cite[Lemma 2.4]{Hu3}), there exists therefore
$N_\gamma\in\N$ such that $b^{N_\gamma}B_1\subset\Img(f'_\gamma)$.
It follows easily that
\set\begin{equation}\label{eq_estimate-coker-f_gamma}
g^\gamma b^{N_\gamma}B_1\subset\Img(f_\gamma).
\end{equation}
By virtue of lemma \ref{lem_lim-one-special-case}(iv), we
deduce that
\set\begin{equation}\label{eq_estimate-coker-phi}
g^\gamma b^{N_\gamma}\Coker\,\phi=0
\end{equation}
where $\phi$ is as in \eqref{eq_very-late}. On the other hand,
since $A_1[1/(gb)]\otimes_{A_1}f_\gamma:A[1/g]^{\oplus n(\gamma)}\to B[1/b]$
is surjective, and since $B[1/b]$ is a projective $A[1/g]$-module,
there exists an $A[1/g]$-linear map $h:B[1/b]\to A[1/g]^{\oplus n(\gamma)}$
such that $f_\gamma\circ h=\one_{B[1/b]}$. Let $\eps_1,\dots,\eps_{n(\gamma)}$
be the canonical basis of the free $A[1/g]$-module
$A[1/g]^{\oplus n(\gamma)}$, and pick $k\in\N$ such that
$h\circ f_\gamma(b^kg^k\eps_i)\in A_1^{\oplus n(\gamma)}$ for every
$i=1,\dots,n(\gamma)$. In view of \eqref{eq_estimate-coker-f_gamma},
we deduce an $A_1$-linear map :
$$
h':B_1\to A_1^{\oplus n(\gamma)}
\qquad
x\mapsto b^{N_\gamma+k}g^{k+\gamma}h(x)
$$
and by construction we have
$f_\gamma\circ h'=b^{N_\gamma+k}g^{k+\gamma}\one_{B_1}$. Then lemma
\ref{lem_lim-one-special-case}(iv) implies that
\set\begin{equation}\label{eq_estimate-ker-phi}
b^{N_\gamma+k}g^{k+\gamma}\Ker\,\phi=0.
\end{equation}
Now, let $e_\bullet:=(e_n~|~n\in\N)\in\lim_{n\in\N}(R^a/b^nR^a)_*$
be the compatible system of diagonal idempotents of the \'etale
$(A_0/b^nA_0)^a$-algebras $B^a_1/b^nB^a_1$, for every $n\in\N$;
by \eqref{eq_estimate-coker-phi} there exists $e'\in R$ such that
$\phi(e')=g^{2\gamma}b^{N_\gamma}e_\bullet$. Then we have
$$
\phi(x\cdot e')=0
\qquad\text{for every $x\in\Ker\,\mu$}.
$$
Combining with \eqref{eq_estimate-ker-phi}, it follows that
$b^{N_\gamma+k}g^{k+\gamma}x\cdot e'=0$ in $R$, for every $x\in\Ker\,\mu$.
But we have already observed that $b$ is regular in $R[1/g]$;
we then easily conclude that the image of $e'$ in $R[1/g]$
agrees with $g^{2\gamma}b^{N_\gamma}e_{B_1/A_0[1/g]}$.

(ii): Pick $n\in\N$ such that $g^ne_{B/A_0[1/g]}$ lies in the
image of the localization $B_1\otimes_{A_1}B_1\to B\otimes_{A_0[1/g]}B$.
Arguing as in the proof of claim \ref{cl_nausee}, we then find
$k\in\N$ and $A_1$-linear maps $B_1\to A_1^{\oplus k}\to B_1$
whose composition is $g^k\cdot\one_{B_1}$. Recall now that
$R^\wedge$ is the limit of the inverse system
$(B_1\otimes_{A_1}B_1/b^iB_1~|~i\in\N)$; then the first
assertion follows from lemma \ref{lem_lim-one-special-case}(iv).

In order to show the second assertion of (ii), let $C_n$
be the \'etale $\sA^a_n$-algebra of finite rank as in the
proof of theorem \ref{th_perf-Abhyankar}, for every $n\in\N$.
Since $\sA_n$ is $b$-adically complete and separated, the
almost projective $\sA^a_n$-module of finite rank 
$S_n:=C_n\otimes_{\sA_n}C_n$ is $b$-adically complete and separated
as well, for every $n\in\N$, by \cite[Claim 5.3.25]{Ga-Ra},
and moreover $\Ann_{S_n}(g)=0$, since $\Ann_{\sA_n}(g)=0$; then
we are easily reduced to checking :

\begin{claim} $(R^\wedge)^a$ represents the limit of the inverse
system of $A_0$-algebras $(S_n~|~n\in\N)$.
\end{claim}
\begin{pfclaim}[] It suffices to show that the
$A^a_0/b^kA^a_0$-algebra $(B_1/b^kB_1\otimes_{A_1/b^kA_1}B_1/b^kB_1)^a$
represents the limit of the induced inverse system
$(C_n/b^kC_n\otimes_{\sA^a_n/b^k\sA^a_n}C_n/b^kC_n~|~n\in\N)$, for every
$k\in\N$. However, the proof of theorem \ref{th_perf-Abhyankar}
shows that the inverse system $(C_n/b^kC_n~|~n\in\N)$ is almost
essentially constant for every $k\in\N$, and the same holds fo
the inverse system $(\sA^a_n/b^k\sA^a_n~|~n\in\N)$, by claim
\ref{cl_ginocchio-quasi-ok}(ii). Then the assertion follows
from proposition \ref{prop_Aretha-RIP}.
\end{pfclaim}
\end{proof}

We conclude with a discussion of some further constructions
that we do not need (though, see remark \ref{rem_alternative}),
but may be useful for other questions, and are related to some
recent work of Y.Andr\'e and others as well.

\sset\subsubsection{}\label{subsec_adjoints-for-PerfTate}
In the situation of \eqref{subsec_almost-perfectoids}, suppose
that $A_0$ is perfectoid; consider the categories :
$$
A_0\tdu\mathsf{cp.Adic}
\qquad
A_0\tdu\Tate
\qquad
A_0\tdu\mathsf{Perf.Adic}
\qquad
A_0\tdu\mathsf{Perf.Tate}
$$
such that :
\begin{itemize}
\item
the objects of $A_0\tdu\mathsf{cp.Adic}$ are the adic ring
homomorphisms $A_0\to B_0$ of complete separated topological
rings, and the morphisms are the maps of $A_0$-algebras
\item
$A_0\tdu\Tate$ is the full subcategory of $A_0\tdu\mathsf{cp.Adic}$
whose objects are the $A_0$-algebras $B_0$ such that the localisation
$B_0\to B_0[b^{-1}]$ is injective
\item
$A_0\tdu\mathsf{Perf.Adic}$ (resp. $A_0\tdu\mathsf{Perf.Tate}$)
is the full subcategory of $A_0\tdu\mathsf{cp.Adic}$ (resp. of
$A_0\tdu\Tate$) whose objects are the perfectoid $A_0$-algebras.
\end{itemize}
Notice that every morphism in these categories is adic. We
wish to exhibit right adjoints for the inclusion functors
\set\begin{equation}\label{eq_from-PerfTate-to-Tate}
A_0\tdu\mathsf{Perf.Adic}\to A_0\tdu\mathsf{cp.Adic}
\qquad
A_0\tdu\mathsf{Perf.Tate}\to A_0\tdu\Tate.
\end{equation}
Thus, consider any $B_0\in\Ob(A_0\tdu\mathsf{cp.Adic})$, set
$\bar A_0:=A_0/pA_0$, $\bar B_0:=B_0/pB_0$, and let
$\bar\beta_\bullet:=(\bar\beta_n~|~n\in\N)\in\bE_0:=\bE(\bar A_0)$
be any element such that $\bar\beta_0\in\bar A_0$ is the class of
$b$. Endow $\bE(\bar B_0)$ with its $\bar\beta_\bullet$-adic
topology $\cT_{\bE(\bar B_0)}$; let also $\cT_{\bar B_0}$ be the
$b$-adic topology of $\bar B_0$. Then the topology of
$\bE(\bar B_0,\cT_{\bar B_0})$ is complete, separated, adic,
and coarser than $\cT_{\bE(\bar B_0)}$ (remark
\ref{rem_topology-of-E}(ii)). By lemma \ref{lem_fontaine}, it
follows that $\cT_{\bE(\bar B_0)}$ is complete and separated.
Then the $A_0$-algebra
$$
B_0^\natural:=A_0\otimes_{W(\bE_0)}W(\bE(\bar B_0),\cT_{\bE(\bar B_0)})
$$
is a perfectoid ring for its $b$-adic topology (example
\ref{ex_perfectoid}(ii)). The structure map $A_0\to B_0$ of
$B_0$, and the map $u_{B_0}:W(\bE(\bar B_0))\to B_0$ induce an
adic map
$$
\eps_{B_0}:B_0^\natural\to B_0
\qquad
1\otimes(x_n~|~n\in\N)\mapsto\sum_{n\in\N}p^n\cdot\bar u_A(x_n^{p^{-n}})
$$
and notice that $\eps_{B_0}$ is an isomorphism if $B_0$ is
perfectoid. Next, let $\beta_\bullet:=(\beta_n~|~n\in\N)\in\bE(A_0)$
be the preimage of $\bar\beta_\bullet$ under the isomorphism
$\bE(A_0)\isom\bE_0$ induced by the projection $A_0\to\bar A_0$.
Hence $\beta_0-b\in pA_0$, so that $\beta_0/b\in A_0^\times$; if
$B_0\in\Ob(A_0\tdu\Tate)$, then the image of $b$ is regular in
$B_0$, so the same holds for the image of $\beta_0$, and finally,
the image of $\beta_\bullet$ is a regular element of $\bE(\bar B_0)$
(proposition \ref{prop_back-and-forth-sec}); by the same token, the
image of $\beta_0$ in $B_0^\natural$ is then regular, and hence the same
holds for the image of $b$. Hence, in this case $B_0^\natural$ is an
object of $A_0\tdu\mathsf{Perf.Tate}$. Moreover, every morphism
$g:B_0\to C_0$ of $A_0\tdu\Tate$ induces an adic map of $A_0$-algebras
$$
g^\natural:B_0^\natural\to C_0^\natural
\qquad
1\otimes(\bar x_n~|~n\in\N)\mapsto
1\otimes(\bE(\bar g)(\bar x_n)~|~n\in\N)
$$
with $\bar g:=g\otimes_\Z\Z/p\Z:\bar B_0\to\bar C_0$. We have
therefore well defined functors
$$
(-)^\natural:A_0\tdu\mathsf{cp.Adic}\to A_0\tdu\mathsf{Perf.Adic}
\qquad
(-)^\natural:A_0\tdu\Tate\to A_0\tdu\mathsf{Perf.Tate}.
$$
Now, let $C\in\Ob(A_0\tdu\Perf)$; to every adic map
$f:C\to B_0^\natural$ let us attach $f^*:=\eps_{B_0}\circ f:C\to B_0$,
and to every adic map $h:C\to B_0$ let us attach
$h^*:=h^\natural\circ\eps^{-1}_C:C\to B_0^\natural$. By a simple
inspection we get $h^{**}=h$. In order to check that $f^{**}=f$,
it suffices to show that $\eps_{B_0^\natural}=\eps_{B_0}^\natural$. Thus,
let $\bar x_\bullet:=(x_n~|~n\in\N)\in W(\bE(\bar B{}_0^\natural))$
be any element; via the natural identification
$\bE(B_0^\natural)\isom\bE(\bar B{}^\natural_0)$, we may view
$\bar x_\bullet$ as a system $x_{\bullet\bullet}:=(x_{n,k}~|~n,k\in\N)$
of elements of $B_0^\natural$ such that $x_{n,k+1}^p=x_{n,k}$ for
every $n,k\in\N$; then we have :
$$
\eps_{B_0^\natural}(1\otimes x_{\bullet\bullet})=\sum_{n\in\N}p^nx_{n,n}
\qquad
\eps^\natural_{B_0}(1\otimes\bar x_\bullet)=
1\otimes(\bE(\bar\eps_{B_0})(\bar x_n)~|~n\in\N).
$$
But we have
$\sum_{n\in\N}p^n\cdot\tau(\bE(\bar\eps_{B_0})(\bar x_n)^{p^{-n}})=
(\bE(\bar\eps_{B_0})(\bar x_n)~|~n\in\N)$ in $W(\bE(B_0^\natural))$,
by virtue of \eqref{eq_new-form}, where $\tau$ denotes the
Teichm\"uller representative. Hence, let
$\bar y_\bullet:=(\bar y_k~|~k\in\N)\in\bE(\bar B{}_0^\natural)$
be any element, which we identify with a unique element
$(y_k~|~k\in\N)\in\bE(B{}_0^\natural)$; we come down to showing
that $1\otimes\tau(\bar y_\bullet)=y_0=\bar u_{B_0^\natural}(y_\bullet)$
in $B_0^\natural$. But we have
$1\otimes\tau(\bar y_\bullet)=u_{B_0^\natural}(\tau(\bar y_\bullet))$,
so the assertion follows from lemma \ref{lem_drop-conditions}(iv).
This shows that the rule $f\mapsto f^*$ yields the sought
adjunctions between \eqref{eq_from-PerfTate-to-Tate} and the
functors $(-)^\natural$.

\sset\subsubsection{}
Keep the notation and assumptions of
\eqref{subsec_adjoints-for-PerfTate}, and consider as well
the categories
$$
(A_0,\fm)^a\tdu\mathsf{Perf.Adic}
\qquad
(A_0,\fm)^a\tdu\mathsf{Perf.Tate}
$$
with the same objects as $A_0\tdu\mathsf{Perf.Adic}$ (resp. as
$A_0\tdu\mathsf{Perf.Tate}$), and whose morphisms $B_0\to C_0$ are
the morphisms $B_0^a\to C_0^a$ of $(A_0,\fm)^a$-algebras. We have
obvious forgetful functors :
\set\begin{equation}\label{eq_carina}
\begin{aligned}
A_0\tdu\mathsf{Perf.Adic}\to & (A_0,\fm)^a\tdu\mathsf{Perf.Adic} \\
A_0\tdu\mathsf{Perf.Tate}\to & (A_0,\fm)^a\tdu\mathsf{Perf.Tate}
\end{aligned}
\qquad
B_0\mapsto B_0^a.
\end{equation}
As explained in \eqref{subsec_make-a-comp-pair-ofalm-str},
we may attach to $\fm$ an ideal $\fm_\bE\subset\bE_0$, and
remark \ref{rem_make-comp-pair-alm-struct}(i) implies that
$(\bE_0,\fm_\bE)$ is a basic setup compatible with $(A_0,\fm)$,
in the sense of definition \ref{def_compat-alm-structures}.
For every $B_0\in\Ob(A_0\tdu\mathsf{cp.Adic})$, we then let $B^a_0$
(resp. $\bE(B_0)^a$) be the $(A_0,\fm)^a$-algebra attached to $B$
(resp. the $(\bE_0,\fm_\bE)^a$-algebras attached to $\bE(B_0)$).

\begin{theorem}\label{th_about-Andre-construction}
Let $B_0\in\Ob(A_0\tdu\mathsf{cp.Adic})$, and endow
$C:=(B_0)^a_{!!}$ with its $b$-adic topology, and
$D:=\bE(B_0)^a_{!!}$ with its $\bar\beta_\bullet$-adic
topology, as in \eqref{subsec_adjoints-for-PerfTate}. We have :
\begin{enumerate}
\item
$(B_0)^a_*\in\Ob(A_0\tdu\mathsf{cp.Adic})$, and if
$B_0\in\Ob(A_0\tdu\Tate)$, then $(B_0)^a_*\in\Ob(A_0\tdu\Tate)$.
\item
There exists a natural isomorphism $\bE(B^a_{0*})\isom(\bE(B_0))^a_*$
of\/ $\bE_0$-algebras.
\item
Suppose moreover that $B_0$ is almost perfectoid; then we have :
\begin{enumerate}
\item
$(B_0)^a_*$ is almost perfectoid.
\item
$\eps^a_{B_0}:(B_0^\natural)^a\to B_0^a$ is an isomorphism.
\item
The completions $C^\wedge$ and $D^\wedge$ of\/ $C$ and $D$ are perfectoid.
\item
$C^\wedge/\Ann_{C^\wedge}(b)\in\Ob(A_0\tdu\mathsf{Perf.Tate})$.
\item
There exists a natural isomorphism $\bE(C^\wedge)\isom D^\wedge$.
\end{enumerate}
\item
The natural morphisms $(C^\wedge)^a\to B^a_0$ and
$(D^\wedge)^a\to\bE(B_0)^a$ are isomorphisms.
\item
If $B_0$ is perfectoid, the same holds for $C^\wedge$ and $D^\wedge$.
\item
The functors \eqref{eq_carina} admit both left and right
adjoints.
\end{enumerate}
\end{theorem}
\begin{proof}(i): If $b\cdot\one_{B_0}$ is an injective map,
the same holds for $b\cdot\one_{B^a_{0*}}$, because the functor
$(-)^a_*$ commutes with limits; the assertion then follows from
the more general :

\begin{claim} Let $(V,\fm_V)$ be a basic setup, $A$ a $V$-algebra,
$I\subset A$ an ideal of finite type, $M$ an $A$-module whose
$I$-adic topology is complete and separated. Then the same holds
for the $I$-adic topology of $M^a_*$.
\end{claim}
\begin{pfclaim} Both functors $(-)^a$ and $(-)_*$ commute
with limits, since they are right adjoints, hence $M^a_*$
is the limit of the inverse system $((M/I^nM)^a_*~|~n\in\N)$.
For every $n\in\N$, let $M_n$ be the image of the natural
map $\pi_n:M^a_*\to(M/I^nM)^a_*$; it follows that $M^a_*$ is
also the limit of the induced inverse system $(M_n~|~n\in\N)$,
and therefore $M^a_*$ is complete and separated for the
linear topology that admits the fundamental system of
open submodules $(\Ker\,\pi_n=(I^nM)^a_*~|~n\in\N)$
(corollary \ref{cor_limits-and-complete}(i)). Since
$I^nM^a_*\subset(I^nM)^a_*$ for every $n\in\N$, the
claim now follows from lemma \ref{lem_fontaine}.
\end{pfclaim} 

(ii): By remark \ref{rem_make-comp-pair-alm-struct}(iii), the
natural map $B_0\to B^a_{0*}$ induces an isomorphism
$\bE(B_0)^a_*\isom\bE(B^a_{0*})^a_*$; in light of (i), we may
then replace $B_0$ by $B^a_{0*}$, and assume from start that
$B_0=B_{0*}^a$, in which case we show that the natural map
$\bE(B_0)\to\bE(B_0)^a_*$ is an isomorphism. To this aim,
consider the system of sets and maps of sets
$S_{\bullet\bullet}:=(S_{ij},\phi^v_{ij},\phi^h_{ij}~|~i,j\in\N)$ with
$S_{ij}:=(B_0/b^jB_0)^a_*$ for every $i,j\in\N$, and where :
\begin{itemize}
\item
$\phi^h_{ij}:S_{i+1,j}\to S_{ij}$ is given by the rule : $x\mapsto x^p$
for every $x\in(B_0/b^jB_0)^a_*$
\item
$\phi^v_{ij}:S_{i,j+1}\to S_{ij}$ is induced by the projection
$B_0/b^{j+1}B_0\to B_0/b^jB_0$.
\end{itemize}
We compute the inverse limit $L$ of the system $S_{\bullet\bullet}$
in two different ways :
$$
L=\lim_{i\in\N}\lim_{j\in\N}S_{ij}\isom\lim_{j\in\N}\lim_{i\in\N}S_{ij}.
$$
Notice that, for fixed $i\in\N$, the system
$S_{i\bullet}:=(S_{ij},\phi^v_{ij}~|~j\in\N)$ consists of rings and
ring homomorphisms, and since the functor $(-)^a_*$ commutes with
limits, the limit of $S_{i\bullet}$ represents $L_i:=B^a_{0*}=B$; then,
the system of maps $(\phi^h_{ij}~|~j\in\N)$ is a natural transformation
$S_{i+1\bullet}\to S_{i\bullet}$ whose limit $\phi^h_i:L_{i+1}\to L_i$ is
obviously the map given by the rule : $x\mapsto x^p$ for every $x\in B$.
Thus $L:=\lim_{i\in\N}(S_{i\bullet},\phi^h_i)$ represents the set
underlying $\bE(B)$. Next, for every $j\in\N$, let $L'_j$ be the
limit of the system of sets $S_{\bullet j}:=(S_{ij},\phi^h_{ij}~|~i\in\N)$;
again, the system $(\phi^v_{ij}~|~i\in\N)$ is a natural
transformation $S_{\bullet j+1}\to S_{\bullet j}$, inducing a map
$\phi^v_j:L'_{j+1}\to L_j$ and $L$ is also $\lim_{j\in\N}(L'_j,\phi^v_j)$.
Notice now that $L'_0$ is the set underlying $\bE(B_0/bB_0)^a_*$,
since the functor $(-)^a_*$ commutes with limits. We are therefore
reduced to checking that $\phi^v_j$ is bijective for every $j\in\N$.
To this aim, it suffices to exhibit a system of maps
$(s_{ij}:S_{i+1,j}\to S_{i,j+1}~|~i\in\N)$ such that
$$
s_{ij}\circ\phi^v_{i+1,j}=\phi^h_{i,j+1}
\qquad
\phi^v_{ij}\circ s_{ij}=\phi^h_{ij}
\qquad
\text{for every $i,j\in\N$}
$$
(details left to the reader). This in turn is achieved by the
following more general :

\begin{claim} Let $(V,\fm_V)$ be a basic setup such that
$\tilde\fm_V:=\fm_V\otimes_V\fm_V$ is a flat $V$-module, $R$
a $V$-algebra, $I\subset R$ an ideal with $pI=0$ and
$x^p=0$ for every $x\in I$. Denote by $\pi:R\to R/I$ the
projection. Then there exists a map $s$ that makes commute
the diagram of sets :
$$
\xymatrix{ R^a_* \ar[r]^-\Phi \ar[d]_{\pi^a_*} &
R^a_* \ar[d]^{\pi^a_*} \\
(R/I)^a_* \ar[r]^-{\bar\Phi} \ar[ru]^-s & (R/I)^a_*
}$$
where $\Phi$ is the map given by the rule : $x\mapsto x^p$
for every $x\in R^a_*$, and likewise for $\bar\Phi$.
\end{claim}
\begin{pfclaim} For every $V$-algebra $S$, denote by
$\sigma_S:R/I\otimes_VS\to R\otimes_VS$ the map defined
as follows. For every $\bar x\in R/I\otimes_VS$, pick
$x\in R\otimes_VS$ such that $\pi(x)=\bar x$, and set
$\sigma_S(x):=x^p$; it is easily seen that this map is
independent of the choice of representatives, and moreover
the rule : $S\mapsto\sigma_S$ yields a homogeneous polynomial
law $\sigma:R/I\leadsto R$ of degre $p$ (see definition
\ref{def_multipolynomial}(ii)). Thus, $\sigma$ factors
through the universal homogeneous polynomial law
$\blambda^p_{R/I}:R/I\leadsto\Gamma^p_V(R/I)$ and a unique
$R$-linear map $\sigma^*:\Gamma^p_V(R/I)\to R$ (see
\eqref{subsec_comultipl-Gamma}). On the other hand, since
$\tilde\fm_V$ is a flat $V$-module, we have a $V$-linear
isomorphism $\omega:\tilde\fm_V\isom\Gamma^p_V(\tilde\fm_V)$
(\cite[(2.1.11)]{Ga-Ra}); then the map given by the rule :
$$
\Hom_V(\tilde\fm_V,R/I)\to\Hom_V(\tilde\fm_V,R)
\qquad
\phi\mapsto(\tilde\fm_V\xrightarrow{\omega}\Gamma^p_V(\tilde\fm_V)
\xrightarrow{\Gamma^p_V(\phi)}\Gamma^p_V(R/I)\xrightarrow{\sigma^*}R)
$$
fulfills the required conditions : details left to the reader.
\end{pfclaim}

(iii.a): In light of (i), it suffices to show that the Frobenius
endomorphism of $R:=B^a_{0*}$ induces an isomorphism
$(R/bR)^a\isom(R/b^pR)^a$; the latter is clear, since the same
holds by assumption for the Frobenius endomorphism of $B_0$, and
the natural map $B_0\to B^a_{0*}$ is an almost isomorphism.

(iii.b): Let $B_1$ be the $p$-integral closure of $B_0$
in $B_0[b^{-1}]$, endow $B_1$ with its $b$-adic topology, and
let $B_1^\wedge$ be the completion of $B_1$. Since $b\cdot\one_{B_1}$
is injective, the same holds for
$(b\cdot\one_{B_1})^\wedge=b\cdot\one_{B_1^\wedge}$ (proposition
\ref{prop_replaces-Mat-Th.8.1}(i)). Then, according
to claim \ref{cl_roads-Portishead} and proposition
\ref{prop_almost-perfectoid}(i), the $A_0$-algebra $B_1^\wedge$
is almost perfectoid and the induced map $i:B_0\to B_1^\wedge$ is
an almost isomorphism. Moreover, $B_1^\wedge$ is $p$-integrally
closed in $B_1^\wedge[b^{-1}]$, by lemma
\ref{lem_crit-p-integr-closed}(ii) and proposition
\ref{prop_complete-f-adic}(iii). We consider the commutative
diagram
$$
\xymatrix{ B_0^\natural \ar[r]^-{\eps_{B_0}} \ar[d]_{i^\natural} &
B_0 \ar[d]^i \\
(B_1^\wedge)^\natural \ar[r]^-{\eps_{B_1^\wedge}} & B_1^\wedge.
}$$
By remark \ref{rem_make-comp-pair-alm-struct}(iii), $i$ induces
an isomorphism $\bE(i)^a:\bE(B_0)^a\to\bE(B^\wedge_1)^a$ of
$(\bE,\fm_\bE)^a$-algebras. Let $\fn\subset A_0$ be the ideal
generated by $\bar u_{A_0}(\fm_\bE)$; it is easily seen that
$(A_0,\fn)$ is a basic setup, and then
$i^\natural:=A_0\otimes_{W(\bE_0)}W(\bE(i))$ is an almost isomorphism
relative to the almost structure furnished by $(A_0,\fn)$. Now,
$i^\natural$ is the inverse limit of the system of maps
$(i_n:=i^\natural\otimes_\Z\Z/p^n~|~n\in\N)$, and in view of
\eqref{eq_same-p-top-closure} we deduce that $i_n$ is an almost
isomorphism also relative to the original almost structure
$(A_0,\fm)$, for every $n\in\N$, and finally, the same holds for
$i^\natural$. Summing up, we are therefore reduced to checking that
$\eps^a_{B_1^\wedge}$ is an isomorphism of $(A_0,\fm)^a$-algebras,
so we may assume from start that $B_0$ is $p$-integrally closed.
Denote by $\gr_\bullet B_1$ the graded ring associated with the
$b$-adic filtration on $B_0$, and likewise define
$\gr_\bullet B_0^\natural$; since both $B_0$ and $B^\natural_0$ are
$b$-adically complete and separated, we are further reduced to
showing that the induced map
$\gr_\bullet\eps_{B_0}:\gr_\bullet B^\natural_0\to\gr_\bullet B_0$ is
an almost isomorphism, and since $b$ is regular in both $B_0$
and $B_0^\natural$, it suffices to prove that the same holds for
$$
\gr_0\eps_{B_0}:B^\natural_0/bB^\natural_0\to B_0/bB_0.
$$

\begin{claim}\label{cl_sharath}
Define $\beta_\bullet\in\bE(A_0)$ as in
\eqref{subsec_adjoints-for-PerfTate}. Then we have :
\begin{enumerate}
\item
The Frobenius endomorphism of $B_0/\beta_nB_0$ induces
an injective map
$$
B_0/\beta_{n+1}B_0\to B_0/\beta_nB_0
\qquad
\text{for every $n\in\N$}.
$$
\item
The kernel of the projection $\bE(B_0)\to B_0/bB_0$ is
$\beta_\bullet\bE(B_0)$.
\item
We have a natural identification
$B^\natural_0/bB^\natural_0\isom\bE(B_0)/\beta_\bullet\bE(B_0)$.
\end{enumerate}
\end{claim}
\begin{pfclaim}(i) follows straightforwardly from the assumption
that $B_0$ is $p$-integrally closed in $B_0[b^{-1}]$, and (ii)
follows easily from (i) : the details are left to the reader.

(iii): Let $\alpha_\bullet\in W(\bE_0)$ be a distinguished element
in the kernel of $u_A:W(\bE_0)\to A_0$; we have $p\in b^p$ by
assumption, and $\alpha_0\in\beta_\bullet^p\bE_0$, by lemma
\ref{lem_before-name}(ii). But the natural isomorphism
$$
\bE(B_0)/\alpha_0\bE(B_0)\isom
\bE_0/\alpha_0\bE_0\otimes_{\bE_0}\bE(B_0)\isom
A_0\otimes_{W(\bE_0)}\bE(B_0)\isom B^\natural_0/pB^\natural_0
$$
maps the class of $\beta_\bullet$ to the class of $b$,
whence the contention.
\end{pfclaim}

By claim \ref{cl_sharath}(ii), the projection
$\bar u_{B_0/pB_0}:\bE(B_0)\isom\bE(B_0/pB_0)\to B_0/bB_0$ induces
an injective map
$$
j:\bE(B_0)/\beta_\bullet\bE(B_0)\to B_0/bB_0
$$
and claim \ref{cl_sharath}(iii) yields a natural identification
of $j$ with $\gr_0\eps_{B_0}$. Thus, it remains only to check that
$\bar u_{B_0/pB_0}^a$ is an epimorphism of $(\bE_0,\fm_\bE)^a$-algebras.
However, recall that $\bE(B_0)$ is the limit of the system of
rings $(R_n~|~n\in\N)$, with $R_n:=B_0/bB_0$ for every $n\in\N$,
and with transition maps given by the Frobenius endomorphism
$\Phi_{B_0/bB_0}$. Thus, $\bar u_{B_0/pB_0}$ is the inverse limit of
the system of maps $(\Phi^n_{B_0/bB_0}:R_n\to B_0/bB_0~|~n\in\N)$,
each of which is almost surjective, since $B_0$ is almost perfectoid
(as usual, we regard $(\Phi^n_{B_0/bB_0})^a$ as a morphism of
$(\bE_0,\fm_\bE)^a$-algebras $R^a_{n,(\Phi^{-n}_{\bE_0})}\to B^a_0/bB^a_0$ :
cp. remark \ref{rem_make-comp-pair-alm-struct}(iii)). To conclude,
it suffices now to invoke lemma \ref{lem_lim-one-special-case}(i).

(iii.c): We consider first the assertion for $C^\wedge$ : set
$\bar A_0:=A_0/bA_0$, and $\bar\fm:=\fm\bar A_0$; by assumption,
the Frobenius endomorphism of $B_0/b^pB_0$ induces an isomorphism
$\bar\Phi{}^a_{B_0}:(B_0/bB_0)^a\isom(B_0/b^pB_0)_{(\Phi)}^a$ of
$(\bar A_0,\bar\fm)^a$-algebras, whence the isomorphism of
$\bar A_0$-algebras
$$
(\bar\Phi_{B_0})^a_{!!}:(B_0/bB_0)_{!!}^a\isom
((B_0/b^pB_0)_{(\Phi)}^a)_{!!}
$$
where $(-)_{!!}$ denotes the left adjoint of the localization
functor $(-)^a:\bar A_0\Alg\to(\bar A_0,\bar\fm)^a\Alg$. On the
other hand, by a direct inspection of the construction of these
left adjoints, one gets natural identifications
$$
C/bC\isom(B_0/bB_0)_{!!}^a
\qquad
C/b^pC\isom((B_0/b^pB_0)_{(\Phi)}^a)_{!!}
$$
where, however, the construction of $C:=(B_0)^a_{!!}$ refers to
the left adjoint of the corresponding localization
$(-)^a:A_0\Alg\to(A_0,\fm)^a\Alg$. It is then easily seen that,
under these identifications, the map $(\bar\Phi_{B_0})^a_{!!}$
corresponds to the ring homomorphism $\bar\Phi_C:C/bC\to C/b^pC$
induced by the Frobenius endomorphism of $C/b^pC$. The latter
is therefore an isomorphism, and so the same holds for the
corresponding ring homomorphism
$\bar\Phi_{C^\wedge}:C^\wedge/bC^\wedge\to C^\wedge/b^pC^\wedge$.
Then, according to corollary \ref{cor_variant-BDS}, in order to
check that $C^\wedge$ is perfectoid, it remains only to show that
$\Ann_{C^\wedge}(b^p)=\Ann_{C^\wedge}(b^{p-1})$, and in light of remark
\ref{rem_Ann-and-completion}(iii), we are reduced to proving that
$\Ann_C(b^p)=\Ann_C(b^{p-1})$. Let $\mu:\tilde\fm\to A_0$ be the
$A_0$-linear map such that $\mu(x\otimes y):=xy$ for every
$x,y\in\fm$; we notice :

\begin{claim}\label{cl_alas}
$\Ann_C(b^n)$ is a quotient of $\Ker(\tilde\fm/b^n\tilde\fm
\xrightarrow{A_0/b^nA_0\otimes_{A_0}\mu}A_0/b^nA_0)$, for all $n\in\N$.
\end{claim}
\begin{pfclaim} We consider the commutative diagram of
$A_0$-modules with exact rows :
$$
\xymatrix{ 0 \ar[r] & \fn \ar[r] \ar[d] &
A_0\oplus(\tilde\fm\otimes_{A_0}B_0) \ar[r] \ar[d] &
(B_0)^a_{!!} \ar[r] \ar[d] & 0 \\
0 \ar[r] & \fn \ar[r] &
A_0\oplus(\tilde\fm\otimes_{A_0}B_0) \ar[r] &
(B_0)^a_{!!} \ar[r] & 0
}$$
whose vertical arrows are scalar multiplication by $b^n$, and where
$\fn$ denotes the image of the map
$\tilde\fm\to A_0\oplus(\tilde\fm\otimes_{A_0}B_0)$ given by the
rule : $z\mapsto(\mu(z),z\otimes 1)$ for every $z\in\tilde\fm$.
By assumption, $b$ is a regular element in both $A_0$ and $B_0$, and
$\tilde\fm$ is a flat $A_0$-module, by \cite[Prop.2.1.7(i)]{Ga-Ra}.
Hence, the central vertical arrow is injective, and the assertion
follows easily from the snake lemma.
\end{pfclaim}

In view of claim \ref{cl_alas}, we are reduced to showing that
$b\cdot\Ker(A_0/b^nA_0\otimes_{A_0}\mu)=0$ for every $n\in\N$.
To this aim, write $\fm=\bigcup_{n\in\N}A_0a_i$ for a system
$(a_i~|~i\in\N)$ of elements of $A_0$ such that for every
$i\in\N$ we have $a_i=c_ia^p_{i+1}$ for some $c_i\in A_0$ (see
\eqref{subsec_make-a-comp-pair-ofalm-str}); by remark
\ref{rem_make-comp-pair-alm-struct}(ii), the $A_0$-module
$\tilde\fm$ is then isomorphic to the inductive limit $L$
of a system $(F_i~|~i\in\N)$ of free rank one $A_0$-modules
$F_i=A_0e_i$, with transition maps $f_i:F_i\to F_{i+1}$ such that
$f_i(e_i)=c_ia_{i+1}^{p-1}e_{i+1}$ for every $i\in\N$. Hence $L$
is the set of equivalence classes $[x,i]$ of pairs $(x,i)$ with
$i\in\N$ and $x\in F_i$, and with this notation, $\mu$ corresponds
to the $A_0$-linear map $L\to A_0$ given by the rule :
$[ye_i,i]\mapsto ya_i$ for every $i\in\N$ and $y\in A_0$.
Hence, let $[ye_i,i]\in L$ and suppose that $ya_i\in b^nA_0$;
we need to show that $[ye_i,i]\in b^{n-1}L$. Pick $k\in\N$ with
$p^k\geq n$; by a simple induction, we see that :
$$
[ye_i,i]=y_k[e_{i+k},i+k]
\qquad\text{with}\qquad
y_k:=yc_ic_{i+1}^p\cdots c_{i+k-1}^{p^{k-1}}a^{p^k-1}_{i+k}.
$$
We are then reduced to checking that $y_k\in b^{n-1}A_0$, and
by theorem \ref{th_criterium-perfect}, the latter is equivalent
to the condition : $y_k^{p^k}\in b^{(n-1)p^k}A_0$. But we have :
$$
y_k^{p^k}=y^{p^k}c_ic_{i+1}^p\cdots c_{i+k-1}^{p^{k-1}}a_i^{p^k-1}=
yc_ic_{i+1}^p\cdots c_{i+k-1}^{p^{k-1}}\cdot(ya_i)^{p^k-1}\in
b^{n(p^k-1)}A_0\subset b^{(n-1)p^k}A_0.
$$
Next, we consider the assertion concerning $D^\wedge$ : the discussion
of \eqref{subsec_adjoints-for-PerfTate} already shows that $\bE(B_0)$
is a perfectoid ring, when endowed with its $ \bar\beta_\bullet$-adic
topology, and moreover the image of $\bar\beta_\bullet$ is regular in
$\bE(B_0)$. Hence, the assertion for $D^\wedge$ follows from the
foregoing case, after replacing the basic setup $(A_0,\fm)$ by
$(\bE,\fm_\bE)$.

(iii.d): Arguing as in the proof of corollary \ref{cor_variant-BDS},
we may assume that $b=b_0$ for some $(b_n~|~n\in\N)\in\bE(A_0)$;
notice also that $C^\wedge$ is reduced, due to (iii.c) and
corollary \ref{cor_perf-are-reduced}(i), hence
$\Ann_{C^\wedge}(b)=\bigcup_{n\in\N}\Ann_{C^\wedge}(b^n)$, so scalar
multiplication by $b$ is an injective endomorphism of
$C^\wedge/\Ann_{C^\wedge}(b)$. Then the assertion follows from
corollary \ref{cor_BDS}.

(iv): The counit of adjunction for the pair of functors
$((-)_{!!},(-)^a)$ is an almost isomorphism $\eps_{B_0}:C\to B_0$,
hence the same for its $b$-adic completion (cp. the proof of
claim \ref{cl_roads-Portishead}). Likewise, recall that $\bE(B_0)$
is $\bar\beta_\bullet$-adically complete and separated (see
\eqref{subsec_adjoints-for-PerfTate}); then the same argument
shows that the counit of adjunction $\eps_{\bE(B_0)}:D\to\bE(B_0)$
induces, after taking $\bar\beta_\bullet$-adic completions, an
almost isomorphism $D^\wedge\to\bE(B_0)$.

(v): We may assume that $b=\beta_0$, where $\beta_\bullet\in\bE(A_0)$
is as in \eqref{subsec_adjoints-for-PerfTate}; also, arguing as
in the proof of (iii.c), we see that $\bar\Phi_{C^\wedge}$ is an
isomorphism, and then claim \ref{cl_better-estimate} and remark
\ref{rem_Ann-and-completion}(iii) reduce to showing that
$\Ann_C(b^p)=\Ann_C(b^p/\beta_1)=\Ann_C(\beta_1^{p^2-1})$. But,
arguing as in the proof of claim \ref{cl_alas} we obtain, for
every $n\in\N$, an exact sequence of $A_0$-modules
$$
\tilde\fm\otimes_{A_0}\Ann_{B_0}(b^n)\to\Ann_C(b^n)\to M
$$
where $M$ is a quotient of $\Ker(A_0/b^nA_0\otimes_{A_0}\mu)$.
The proof of (iii.c) then shows that $bM=0$, and on the other
hand $\beta_1\cdot\Ann_{B_0}(b^n)=0$, by corollary \ref{cor_BDS};
thus $\beta_1^{p+1}\cdot\Ann_C(b^n)=0$, and since $p^2-1\geq p+1$,
the assertion follows.

(vi): Let $X\in\Ob(A_0\tdu\mathsf{Perf.Tate})$ and
$Y\in\Ob((A_0,\fm)^a\tdu\mathsf{Perf.Tate})$, and set
$Z:=(Y^a_{!!})^\wedge$; we know already from (i) and (iii.d)
that $Y^a_*\in\Ob(A_0\tdu\Tate)$ and
$Z/\Ann_Z(b)\in\Ob(A_0\tdu\mathsf{Perf.Tate})$. Then, in view
of the discussion of \eqref{subsec_adjoints-for-PerfTate}
we get natural bijections
$$
\begin{aligned}
\Hom_{(A_0,\fm)^a\tdu\mathsf{Perf.Tate}}(X^a,Y^a)\isom
\Hom_{A_0\tdu\Tate}(X,Y^a_*)\isom
\Hom_{A_0\tdu\mathsf{Perf.Tate}}(X,(Y^a_*)^\natural) \\
\Hom_{(A_0,\fm)^a\tdu\mathsf{Perf.Tate}}(Y^a,X^a)\isom
\Hom_{A_0\tdu\Alg}(Y^a_{!!},X)\isom
\Hom_{A_0\tdu\mathsf{Perf.Tate}}(Z/\Ann_Z(b),X)
\end{aligned}
$$
so the rule $Y\mapsto(Y^a_*)^\natural$ (resp. $Y\mapsto Z/\Ann_Z(b)$)
yields the sought right (resp. left) adjoint for the functor
$A_0\tdu\mathsf{Perf.Tate}\to(A_0,\fm)^a\tdu\mathsf{Perf.Tate}$.
The same argument shows that the rule $Y\mapsto(Y^a_*)^\natural$
also provides a right adjoint for the forgetful functor
$A_0\tdu\mathsf{Perf.Adic}\to(A_0,\fm)^a\tdu\mathsf{Perf.Adic}$,
and using (v) we easily see that the rule $Y\mapsto Z$ yields a
left adjoint for the same functor.

(iii.e): To begin with, let us notice :

\begin{claim}\label{cl_gas-shoes}
The functors \eqref{eq_carina} induce equivalences between
$(A_0,\fm)^a\tdu\mathsf{Perf.Adic}$ (resp.
$(A_0,\fm)^a\tdu\mathsf{Perf.Tate}$) and the localization of
$A_0\tdu\mathsf{Perf.Adic}$ (resp. $A_0\tdu\mathsf{Perf.Tate}$)
obtained by inverting all almost isomorphisms.
\end{claim}
\begin{pfclaim} Recall that $\eps_{B_0}:B_0^\natural\to B_0$
is an isomorphism if $B_0$ is perfectoid (see
\eqref{subsec_adjoints-for-PerfTate}); we deduce for every
$Y\in\Ob(A_0\tdu\mathsf{Perf.Adic})$ an isomorphism
$((Y^a_*)^\natural)^a\isom(Y^a_*)^a\isom Y^a$ in the category
$(A_0,\fm)^a\tdu\mathsf{Perf.Adic}$. Combining with (vi) and
proposition \ref{prop_fullfaith-adjts}(iii), we deduce that
the right adjoints to both of the forgetful functors
\eqref{eq_carina} are fully faithful. Then the claim follows
from proposition \ref{prop_local-and-adj-pairs}.
\end{pfclaim}

From (vi) and remark \ref{rem_make-comp-pair-alm-struct}(iii)
we get a diagram of categories :
$$
\xymatrix@C+10pt{ A_0\tdu\mathsf{Perf.Adic} \ar[r]^-{\bE(-)}
\ar@<.5ex>[d]^{(-)^a} &
\bE_0\tdu\mathsf{Perf.Adic} \ar@<-.5ex>[d]_{(-)^a} \\
(A_0,\fm)^a\tdu\mathsf{Perf.Adic} \ar[r]^-{\bE(-)} \ar@<.5ex>[u] &
(\bE_0,\fm_\bE)^a\tdu\mathsf{Perf.Adic} \ar@<-.5ex>[u]
}$$
whose upward arrows are left adjoints to the downward arrows,
and whose top horizontal arrow is an equivalence; in light of
claim \ref{cl_gas-shoes}, it follows easily that the bottom
horizontal arrow is an equivalence as well. Clearly the
horizontal arrows commute with the downward arrows; then the
horizontal arrows also commute with the upward arrows, up to
isomorphism of functors (example \ref{ex_presheaves-on-sets}(ii)).
Hence, if $R\in\Ob(A_0\tdu\mathsf{Perf.Adic})$, we have a natural
ring isomorphism
$$
\omega_R:\bE((R^a_{!!})^\wedge)\isom(\bE(R)^a_{!!})^\wedge.
$$
Next, if $B_0$ is almost perfectoid, from (iv) and remark
\ref{rem_make-comp-pair-alm-struct}(iii) we deduce an almost
isomorphism $\bE((B_0)^a_{!!})^\wedge)\to\bE(B_0)$, whence a ring
isomorphism
$$
\tau_{B_0}:(\bE((B^a_{0!!})^\wedge)^a_{!!})^\wedge\isom(\bE(B_0)^a_{!!})^\wedge.
$$
Lastly, $R:=(B^a_{0!!})^\wedge$ is perfectoid, by (v), whence
the isomorphism
$$
\tau_{B_0}\circ\omega_R:\bE((R^a_{!!})^\wedge)\isom
(\bE(B_0)^a_{!!})^\wedge.
$$
Thus, we are reduced to checking that the natural map
$(((B^a_{0!!})^\wedge)^a_{!!})^\wedge\to(B^a_{0!!})^\wedge$ is an
isomorphism. This in turn will follow, once we have shown
that the natural map $((B^a_{0!!})^\wedge)^a_{!!}\to B^a_{0!!}$ is
an isomorphism. But the latter is a straightforward consequence
of (iv).
\end{proof}

\section{Applications}\label{chap_applications}
\subsection{Model algebras}\label{sec_model-alg}
We let $(K,|\cdot|)$ be a valued field of characteristic $0$
with value group $\Gamma_{\!K}$ of rank one, such that the
residue field $\kappa$ of $K^+$ has characteristic $p>0$.

\begin{definition}
The category $\sMA_K$ of {\em model $K^+$-algebras\/} consists
of all the pairs $(A,\Gamma)$, where $(\Gamma,+)$ is an integral
monoid, and $A$ is a $\Gamma$-graded $K^+$-algebra $A$, fulfilling
the following conditions :
\begin{enumerate}
\alphaenu
\item[(MA1)]
$\gr_\gamma A$ is a torsion-free $K^+$-module, with
$\dim_K\gr_\gamma A\otimes_{K^+}K=1$, for every $\gamma\in\Gamma$.
\item[(MA2)]
$\gr_\alpha A\cdot\gr_\beta A\neq 0$ for every
$\alpha,\beta\in\Gamma$, and $(\gr_\gamma A)^n=\gr_{n\gamma}A$
for every $n\in\N$ and $\gamma\in\Gamma$.
\item[(MA3)]
$\gr_0A=K^+$.
\item[(MA4)]
$\Gamma$ is saturated and $\Gamma^\gp$ is a $\Z[1/p]$-module.
\end{enumerate}
The morphisms $(A,\Gamma)\to(A',\Gamma')$ in $\sMA_K$ are the
pairs $(f,\phi)$, where $\phi:\Gamma\to\Gamma'$ is a morphism
of monoids, and $f:A\to\Gamma\times_{\Gamma'}A'$ is a morphism
of $\Gamma$-graded $K^+$-algebras (see \eqref{subsec_monoid-graded}).
\end{definition}

\begin{example} Let $M$ be an integral monoid, $N\to M$ an exact
and injective morphism of monoids, $N\to K^+\!\setminus\!\{0\}$
a morphism of monoids, and suppose that :
\begin{itemize}
\item
$M$ is {\em divisible}, {\em i.e.} the $k$-Frobenius endomorphism
of $M$ is surjective for every $k>0$
\item
$M^\gp/N^\gp$ is a $\Q$-vector space.
\end{itemize}
Let $(\Gamma,+)$ be the image of $M$ in $M^\gp/N^\gp$, and denote
by $I$ the nilradical of the $K^+$-algebra $A:=M\otimes_NK^+$;
then $A$ is a $\Gamma$-graded $K^+$-algebra, and $I$ is a
$\Gamma$-graded ideal (proposition \ref{prop_notnilpo}(ii)). We
claim that $I\otimes_{K^+}K=0$ and $(A/I,\Gamma)$ is a model
$K^+$-algebra. Indeed, (MA4) is immediate, and it is easily seen
that $\gr_{n\gamma}A=(\gr_\gamma A)^n$ for every $\gamma\in\Gamma$
and $n\in\N$. Moreover, since the map $N\to M$ is exact, the
kernel of the map $M\to M^\gp/N^\gp$ equals $N$, so $(A/I)_0=K^+$,
{\em i.e.} (MA3) holds as well. The remaining assertions can be
checked after tensoring with $K$ : namely, we have to show that
the $\Gamma$-graded $K$-algebra $A_K:=M\otimes_NK$ is reduced,
with $\dim_K\gr_\gamma A_K=1$ for every $\gamma\in\Gamma$, and
$\gr_\alpha A_K\cdot\gr_\beta A_K\neq 0$ for every
$\alpha,\beta\in\Gamma$. However, the morphism $N\to K^+$ extends
uniquely to a group homomorphism $N^\gp\to K^\times$, and
$A_K=(N^{-1}M)\otimes_{N^\gp}K$; hence, we may assume that $N$ is a
group, in which case $\Gamma$ is the set-theoretic quotient of $M$
by the translation action of $N$ (lemma \ref{lem_special-p-out}(iii)).
In this case, choose a representative $\gamma^*\in M$ for every
$\gamma\in\Gamma$; it follows easily that
$(\gamma^*\otimes 1~|~\gamma\in\Gamma)$ is a basis of the $K$-vector
space $A_K$, whence (MA1), and
$(\alpha^*\otimes 1)\cdot(\beta^*\otimes 1)$ is a non-zero multiple
of $(\alpha+\beta)^*\otimes 1$, which yields (MA2). Lastly, since
the nilradical $I_K$ of $A_K$ is $\Gamma$-graded (proposition
\ref{prop_notnilpo}(ii)), we also deduce that $I_K=0$, as sought.
\end{example}

\begin{remark}\label{rem_tensor-prod-in-A}
(i)\ \
The category $\sMA_K$ admits a tensor product, defined by the rule :
$$
(A,\Gamma)\otimes(A',\Gamma'):=(A\otimes_{K^+}A',\Gamma\oplus\Gamma').
$$

(ii)\ \
Let $(A,\Gamma)$ be any model $K^+$-algebra, and suppose that
$\Gamma=\Gamma_{\!1}\oplus\Gamma_{\!2}$ is a given decomposition
of $\Gamma$ as direct sum of monoids. Then it is easily
seen that $\Gamma_{\!1}$ and $\Gamma_{\!2}$ are integral and
saturated, and they fulfill axiom (MA4). There follows a morphism
of model $K^+$-algebras :
\set\begin{equation}\label{eq_decompose-in-A}
(\Gamma_{\!1}\times_\Gamma A,\Gamma_{\!1})\otimes
(\Gamma_{\!2}\times_\Gamma A,\Gamma_{\!2})\to(A,\Gamma)
\end{equation}

(iii)\ \
In the situation of (ii), set $A_K:=A\otimes_{K^+}K$. Then it
is easily seen that \eqref{eq_decompose-in-A} induces an
isomorphism of $K$-algebras
$$
(\Gamma_{\!1}\times_\Gamma A_K)\otimes_K(\Gamma_{\!2}\times_\Gamma A_K)
\isom A_K.
$$
In general, \eqref{eq_decompose-in-A} need not be an isomorphism.
However, \eqref{eq_decompose-in-A} is an isomorphism, if either
$\Gamma_{\!1}$ or $\Gamma_{\!2}$ is a torsion abelian group.
Indeed, in any case the induced maps
$$
\gr_\alpha A\otimes_{K^+}\gr_\beta A\to\gr_{\alpha+\beta}A
\qquad
a_1\otimes a_2\mapsto a_1a_2
$$
will be injective for every $\alpha\in\Gamma_{\!1}$ and
$\beta\in\Gamma_{\!2}$. To check surjectivity, let
$a\in\gr_{\alpha+\beta}A$ be any element, and say that $n\beta=0$;
then $a^n\in\gr_{n\alpha}A$, and by (MA2) we know that there
exist $a_1\in\gr_\alpha A$ and $x\in K^+$ such that
$a^n=a_1^nx$. It follows that $a=a_1a_2$ for some
$a_2\in\gr_\beta A_K$ such that $a_2^n\in K^+$. Then
$a_2\in\gr_\beta A$, whence the claim.

(iv)\ \
Let $(E,|\cdot|_E)$ be any valued field extension of $(K,|\cdot|)$,
such that $|\cdot|_E$ is also a valuation of rank one. Then
we have an obvious base change functor :
$$
\sMA_K\to\sMA_E
\quad : \quad
(A,\Gamma)\mapsto(A\otimes_{K^+}E^+,\Gamma).
$$
\end{remark}

\sset\subsubsection{}\label{subsec_not-so-bad}
Let $(A,\Gamma)$ be any model $K^+$-algebra. For any subset
$\Delta\subset\Gamma$, set
$$
A_\Delta:=\Delta\times_\Gamma A
\qquad
A_{\Delta,K}:=A_\Delta\otimes_{K^+}K
\qquad
\Delta[1/p]:=\bigcup_{n\in\N}\{\gamma\in\Gamma~|~p^n\gamma\in\Delta\}.
$$
Clearly, if $\Delta$ is a submonoid of $\Gamma$, then $A_\Delta$
is a $K^+$-subalgebra of $A$, and $(A_\Delta,\Delta)$ is a
subobject of $(A,\Gamma)$, provided $\Delta$ satisfies (MA4).
On the other hand, if $\Delta$ is saturated, it is easily seen
that the same holds for $\Delta[1/p]$, and then the latter does
satisfy (MA4).

Also, if $\Delta$ is an ideal of $\Gamma$, then clearly
$A_\Delta$ is an ideal of $A$.

Set $K^*:=K^+\setminus\!\{0\}$ and let
$A^*_\gamma:=\gr_\gamma A\!\setminus\!\{0\}$ for every
$\gamma\in\Gamma$; for any submonoid $\Delta\subset\Gamma$,
we deduce a sequence of morphisms of integral monoids :
\set\begin{equation}\label{eq_exact-mons}
1\to K^*\xrightarrow{\ \phi_\Delta\ }A_\Delta^*\to\Delta\to 1
\qquad
\text{where $A^*_\Delta:=\bigoplus_{\gamma\in\Delta}A^*_\gamma$}
\end{equation}
(and where the direct sum is formed in the category of $K^*$-modules),
such that $\eqref{eq_exact-mons}^\gp$ is a short exact complex
of abelian groups; then (MA2) implies that $\phi_\Delta$ is saturated
(proposition \ref{prop_crit-saturated}). Also, (MA1) implies that
each $A^*_\gamma$ is a filtered union of free cyclic $K^*$-modules,
hence $\phi_\Delta$ is also integral. Furthermore, (MA1) and (MA2)
imply that $A_\Delta^{*\gp}$ is $\Delta^\gp$-graded, and
\set\begin{equation}\label{eq_A-group}
\gr_\gamma A^{*\gp}_\Delta=A_\gamma^*\otimes_{K^*}K^\times
\qquad
\text{for every $\gamma\in\Delta$}.
\end{equation}
(details left to the reader). We shall just write $A^*$
and $A_K$ instead of $A^*_\Gamma$, and respectively $A_{\Gamma,K}$.

\sset\subsubsection{}\label{subsec_finally-purity}
Let now $(A,\Gamma)$ be any model $K^+$-algebra; set
$S:=\Spec\,K^+$, $X:=\Spec\,A$, and let $\bar x$ be
any geometric point of $X$, localized on the closed
subset $Z:=X\times_S\Spec\,\kappa$. Let $\cT_{A,p}$ be
the $p$-adic topology on $A$. Suppose furthermore, that
$K^+$ is deeply ramified (see \cite[Def.6.6.1]{Ga-Ra}),
and let $(K^+,\fm_K)$ be the standard setup associated
with $K^+$ (see \cite[\S6.1.15]{Ga-Ra}); then we have
the corresponding sheaf $\cO_{\!X(\bar x)}^a$ of
$K^{+a}$-algebras on $X(\bar x)$, and we may state the
following {\em almost purity\/} theorem :

\begin{theorem}\label{th_finally-purity}
With the notation of \eqref{subsec_finally-purity}, the
following holds :
\begin{enumerate}
\item
$(A,\cT_{A,p})$ is a formal perfectoid ring.
\item
The pair $(X(\bar x),Z(\bar x))$ is almost pure.
\end{enumerate}
\end{theorem}
\begin{proof}(i): Let $\pi\in K^+$ be any element such that
$\pi^pK^+=pK^+$; denote by $A^\wedge$ the completion of
$(A,\cT_{A,p})$, and set $I:=\pi A^\wedge$. Then the topology
of $A^\wedge$ agrees with the $I$-adic topology (remark
\ref{rem_completion-of-topring}(iv)), and the image of
$\pi$ is a regular element in $A^\wedge$, by proposition
\ref{prop_replaces-Mat-Th.8.1}(i) (the details shall be
left to the reader). Moreover we have $pA^\wedge=I^{(p)}$
(notation of definition \ref{def_perfectoid}); by theorem
\ref{th_regular-seq-criterion}, it then suffices to check
that the Frobenius endomorphism $\Phi_{A/pA}$ of
$A^\wedge/pA^\wedge=A/pA$ induces an isomorphism
$A/\pi A\isom A/pA$. We are then further reduced to showing
that $\Phi_{A/pA}$ restricts to a bijection
$$
\gr_\gamma A/\pi\cdot\gr_\gamma A\isom
\gr_{p\gamma}A/p\cdot\gr_{p\gamma}A
\qquad
\text{for every $\gamma\in\Gamma$}.
$$
However, since $\Gamma_{\!K}$ is of rank one, there exists
a sequence $(a_n~|~n\in\N)$ of elements of $\gr_\gamma A$
such that $K^+a_n\subset K^+a_{n+1}$ for every $n\in\N$
and $\gr_\gamma A=\bigcup_{n\in\N}K^+a_n$. Then
$\gr_{p\gamma}A=\bigcup_{n\in\N}K^+a_n^p$. By
\cite[Prop.6.6.6]{Ga-Ra} the Frobenius endomorphism of $K^+$
induces an isomorphism $K^+/\pi K^+\isom K^+/pK^+$. It follows
that $\Phi_{A/pA}$ induces a bijection
$(K^+a_n)/\pi(K^+a_n)\isom(K^+a_n^p)/p(K^+a_n^p)$ for every
$n\in\N$, whence the contention.

(ii) follows from (i), and theorems \ref{th_formal-perf}(i)
and \ref{th-alm-purity-form-perfectoid}.
\end{proof}

\sset\subsubsection{}\label{subsec_SS}
The inclusion map $A^*\to A$ is a morphism of (multiplicative)
monoids, hence induces a log structure on $X:=\Spec\,A$ (see
\eqref{subsec_Konstant}). We shall denote
$$
\SS(A,\Gamma):=(X,(A^*_X)^{\log})
$$
the resulting log scheme. Clearly, every morphism
$(f,\phi):(A,\Gamma)\to(A',\Gamma')$ of model algebras
induces a morphism of log schemes
$$
\SS(f,\phi):\SS(A',\Gamma')\to\SS(A,\Gamma).
$$
Especially, by lemma \ref{lem_integr-flat}(iv), the map
$\phi_\Gamma$ yields a saturated morphism
\set\begin{equation}\label{eq_ss-aturated}
\SS(A,\Gamma)\to\SS(K^+):=\SS(K^+,\{1\}).
\end{equation}
Also, if $(K,|\cdot|)\to(E,|\cdot|_E)$ is an extension of
rank one valued fields, the inclusion
$A^*\subset(A\otimes_{K^+}E^+)^*$ induces a morphism of log
schemes
$$
\SS(A\otimes_{K^+}E^+,\Gamma)\to\SS(A,\Gamma)
$$
for any model $K^+$-algebra $(A,\Gamma)$, and we remark that
the resulting diagram of log schemes
$$
\xymatrix{
\SS(A\otimes_{K^+}E^+,\Gamma) \ar[r] \ar[d] & \SS(A,\Gamma) \ar[d] \\
\SS(E^+) \ar[r] & \SS(K^+) 
}$$
is cartesian. Indeed, since \eqref{eq_another-left-adj} is
right exact, it suffices to check that the natural map
$$
A^*\otimes_{K^*}E^*\to(A\otimes_{K^+}E^+)^*
$$
is an isomorphism, which is clear.

\begin{remark}\label{rem_good-luck-wt}
(i)\ \
Let $(A,\Gamma)$ be any model $K^+$-algebra, and suppose that
$\Delta_0\subset\Gamma$ is a fine and saturated submonoid,
such that $\Delta_0^\gp$ is torsion-free. In this case,
$\Delta_0^\gp$ is a free abelian group of finite rank, hence
$\eqref{eq_exact-mons}^\gp$ admits a splitting
$\sigma:\Delta_0^\gp\to A^{*\gp}$. Then, using
\eqref{eq_A-group} and (MA1), it is easily seen that the rule
$\gamma\mapsto\sigma(\gamma)$ extends to an isomorphism of
$\Delta_0$-graded $K$-algebras
$$
K[\Delta_0]\isom A_{\Delta_0,K}.
$$

(ii)\ \
In the situation of (i), suppose furthermore, that $K$ is
algebraically closed. Pick any $x\in K^\times$ such that
$\gamma:=|x|\neq 1$, and let $\La\gamma\Ra\subset\Gamma_{\!K}$
be the subgroup generated by $\gamma$; we define a group
homomorphism $\La\gamma\Ra\to K^\times$ by the rule :
$\gamma^k\mapsto x^k$ for every $k\in\Z$. Since $K^\times$
is divisible, the latter map extends to a group homomorphism
\set\begin{equation}\label{eq_ext-to-gamma_K}
\Gamma_{\!K}\to K^\times
\end{equation}
and since $\Gamma_{\!K}$ is a group of rank one, it is easily
seen that \eqref{eq_ext-to-gamma_K} is a right inverse for
the valuation map $|\cdot|:K^\times\to\Gamma_{\!K}$, whence
a decomposition :
$$
K^\times\isom(K^+)^\times\oplus\Gamma_{\!K}.
$$
On the other hand, set
$A^*_{\Delta_0,K}:=A^*_{\Delta_0}\otimes_{K^*}K^\times$;
from (i) we deduce an isomorphism of $\Delta_0$-graded monoids :
$$
A^*_{\Delta_0,K}\isom\Delta_0\oplus K^\times.
$$
Combining these two isomorphisms, we deduce a surjection
$\tau:A^*_{\Delta_0,K}\to(K^+)^\times$ which is a left inverse to
the inclusion $(K^+)^\times\to A^*_0$. Then, for every
$\gamma\in\Delta_0$, let us set $C_\gamma:=A_\gamma^*\cap\Ker\,\tau$;
there follows a (non-canonical) isomorphism of $\Delta$-graded monoids :
$$
A^*_{\Delta_0}\isom(K^+)^\times\oplus C
\qquad\text{where}\quad
C:=\bigoplus_{\gamma\in\Delta_0}C_\gamma\subset A^*_{\Delta_0}
\quad\text{and}\quad
C^\gp\simeq\Delta_0^\gp\oplus\Gamma_{\!K}
$$
(details left to the reader).

(iii)\ \
Suppose that $\Gamma^\gp$ is torsion-free, and $K$ is still
algebraically closed. Let us set :
$$
\Delta_n:=\{\gamma\in\Gamma~|~p^n\gamma\in\Delta_0\}
\qquad
\text{for every $n\in\N$}.
$$
It is easily seen that $\Delta_n$ is still fine and saturated,
and since $K^\times$ is divisible, we may extend inductively
the splitting $\sigma$ of (i) to a system of homomorphisms
$$
\sigma_n:\Delta^\gp_n\to A^{*\gp}
\qquad
\text{such that}
\quad
\sigma_{n+1|\Delta^\gp_n}=\sigma_n
\quad
\text{for every $n\in\N$}
$$
whence a compatible system of isomorphisms
\set\begin{equation}\label{eq_oupla}
K[\Delta_n]\isom A_{\Delta_n,K}
\qquad
\text{for every $n\in\N$}.
\end{equation}
Proceeding as in (ii), we deduce a compatible system of
isomorphisms of $\Delta_n$-graded monoids
$$
A^*_{\Delta_n}\isom(K^+)^\times\oplus C^{(n)}
\qquad
\text{such that $C^{(n)}\subset C^{(n+1)}$ for every $n\in\N$.}
$$
In this situation, notice that the $p$-Frobenius automorphism
of $\Gamma^\gp$ restricts to an isomorphism
$$
\Delta_{n+1}\isom\Delta_n
\qquad
\text{for every $n\in\N$}.
$$
Likewise, since $\Gamma_{\!\!K}$ is $p$-divisible, taking $p$-th
powers induces isomorphisms
$$
C^{(n+1)}\isom C^{(n)}
\qquad
\text{for every $n\in\N$.}
$$
\end{remark}

\begin{definition}\label{def_another-fam}
Let $(B,\Delta)$ be model $K^+$-algebra. We say that
$(B,\Delta)$ is {\em small}, if the following conditions hold :
\begin{enumerate}
\alphaenu
\item
$\gr_\gamma B$ is a finitely generated $K^+$-module for
every $\gamma\in\Delta$.
\item
$\Delta=\Delta_0[1/p]$ for some fine and saturated submonoid
$\Delta_0$.
\item
$\Delta_0\times_\Delta B$ is a finitely generated $K^+$-algebra.
\end{enumerate}
\end{definition}

\begin{remark}
Notice that condition (a) of definition \ref{def_another-fam}
and axiom (MA1) imply that $\gr_\gamma B$ is a free $K^+$-module
of rank one, for every small model agebra $(B,\Delta)$ and every
$\gamma\in\Delta$. Moreover, actually conditions (b) and (c)
(together with axioms (MA1) and (MA2)) imply condition (a).
Indeed, (c), (MA1) and proposition \ref{prop_four-year-later}(ii)
imply that $\gr_\gamma B$ is a free $K^+$-module of rank one,
for every $\gamma\in\Delta_0$. Now, in case $n\gamma\in\Delta_0$,
(MA2) implies that $\gr_{n\gamma}B$ is generated by an element
of the form $z=x_1\cdots x_n$, for certain
$x_1,\dots,x_n\in\gr_\gamma B$. Suppose that $x_1=ay$ for some
$a\in K^+$ and $y\in\gr_\gamma B$; then $z$ is divisible by
$a$ in $\gr_{n\gamma}B$, so $a\in(K^+)^\times$, {\em i.e.}
$x_1$ generates $\gr_\gamma B$. In view of (b), for every
$\gamma\in\Delta$ we may find $k\geq 0$ such that
$p^k\gamma\in\Delta_0$, so (a) follows.
\end{remark}

\begin{remark}\label{rem_decompose-in-B}
Let $(B,\Delta)$ be any small model algebra, so that conditions
(b) and (c) of definition \ref{def_another-fam} are satisfied for
some submonoid $\Delta_0\subset\Delta$. 

(i)\ \
Choose a decomposition $\Delta^\times_0=G\oplus H$, where $G$
is a free abelian group, and $H$ is the torsion subgroup of
$\Delta_0$; we have an isomorphism
$$
\Delta_0\simeq\Lambda\oplus H
\qquad
\text{with\ \ $\Lambda:=\Delta_0^\sharp\oplus G$}
$$
(lemma \ref{lem_decomp-sats}), which induces a decomposition :
$$
\Delta\simeq\Lambda[1/p]\oplus H
$$
inducing, in turn, an isomorphism of model $K^+$-algebras :
$$
(B,\Delta)\isom
(B_{\Lambda[1/p]},\Lambda[1/p])\otimes(B_H,H)
$$
(remark \ref{rem_tensor-prod-in-A}(iii)). Notice that
both $B_\Lambda$ and $B_H$ are finitely generated $K^+$-algebras
(proposition \ref{prop_four-year-later}(i)). In other words,
every small model algebra can be written as the tensor product
of two small model algebras $(B',\Delta')$ and $(B'',\Delta'')$,
such that $\Delta'{}^\gp$ is torsion-free, and $\Delta''$
is a finite abelian group whose order is not divisible by $p$.

(ii)\ \
For every $n\in\N$, set
$\Delta_n:=\{\gamma\in\Delta~|~p^n\gamma\in\Delta_0\}$.
Then $B_n:=\Delta_n\times_\Delta B$ is a finitely generated
$K^+$-algebra for every $n\in\N$. Indeed, in view of (i),
it suffices to check the assertion in case $\Delta^\gp$ is
a torsion-free abelian group. Now, pick a system $x_1,\dots,x_k$
of homogeneous generators of the $K^+$-algebra $B_0$; by (MA2),
for every $i=1,\dots,k$ there exist a homogeneous element
$y_i\in B_n$ and $u_i\in(K^+)^\times$ such that $u_iy_i^{p^n}=x_i$.
Let $z\in B_n$ be any homogeneous element; then
$z^{p^n}=vx_1^{t_1}\cdots x_k^{t_k}$ for some $v\in K$ and
$t_1,\dots,t_k\in\N$. Set $y:=y_1^{t_1}\cdots y_k^{t_k}$
and $u:=u_1^{t_1}\cdots u_k^{t_k}\in K^\times$; then
$y^{p^n}uv=z^{p^n}$, and since $\Delta^\gp$ is torsion-free,
we deduce that $yw=z$ for some $w\in K^+$, {\em i.e.} the
system $y_1,\dots,y_k$ generates the $K^+$-algebra $B_n$.

(iii)\ \
Suppose that $\Delta$ is a finite group whose order is not
divisible by $p$. Then we claim that $B$ is a finite \'etale
$K^+$-algebra.
Indeed, in view of remark \ref{rem_tensor-prod-in-A}(iii),
it suffices to verify the assertion for $\Delta$ a cyclic
finite group, say of order $n$, with $(n,p)=1$; in the latter
case, (MA2) implies that $B\simeq K^+[X]/(X^n-u)$ for some
$u\in(K^+)^\times$, whence the contention.

(iv)\ \
Let $F\subset\Delta$ be any face. Then $(B_F,F)$ is a small
model algebra as well. Indeed, by (i), it suffices to consider
the case where $\Delta^\gp$ is a torsion-free abelian group.
In this case, set $F_0:=F\cap\Delta_0$; it is easily seen
that $F=F_0[1/p]$, and $F_0$ is a fine and saturated monoid, by
lemma \ref{lem_face}(ii) and corollary \ref{coro_satur-face}(ii).
Moreover, $B_{F_0}=F_0\times_FB_{\Delta_0}$ is a finitely
generated $K^+$-algebra, by proposition
\ref{prop_four-year-later}(i), whence the contention.

(v)\ \
Suppose moreover, that $K$ is algebraically closed, and
let $(B,\Delta)$ be any small model $K^+$-algebra. Then
we have a (non-canonical) isomorphism :
$$
K[\Delta]\isom B_K.
$$
To exhibit such an isomorphism, we may -- in light of (i) --
assume that $\Delta^\gp$ is either torsion-free, or a finite
group of order not divisible by $p$. In the latter case,
the assertion follows easily from (iii). In case $\Delta^\gp$
is torsion-free, the sought isomorphism is the colimit of
the system of isomorphisms \eqref{eq_oupla}.
\end{remark}

\begin{lemma}\label{lem_normal-norman}
Let $(A,\Gamma)$ be a model $K^+$-algebra, and denote by
$\cF(\Gamma)$ the filtered family of all fine and saturated
submonoids of\/ $\Gamma$. We have :
\begin{enumerate}
\item
$A=\colim_{\Delta\in\cF(\Gamma)}A_\Delta$.
\item
Suppose that $\Gamma^\gp$ is a torsion-free abelian group.
Then $A$ is a normal domain.
\end{enumerate}
\end{lemma}
\begin{proof}(i) is an immediate consequence of corollary
\ref{cor_fragment-Gordon}(ii).

(ii): We show first the following :

\begin{claim}\label{cl_first-after-K}
Suppose that $\Gamma^\gp$ is a torsion-free abelian group.
Then $A_K$ is a normal domain.
\end{claim}
\begin{pfclaim} In view of (i), it suffices to show that
$A_{\Delta,K}$ is a normal domain, when $\Delta\subset\Gamma$
is fine and saturated.
The latter assertion follows from remark \ref{rem_good-luck-wt}(i)
and theorem \ref{th_structure-of-satu}(iii).
\end{pfclaim}

Let $A^\nu$ be the integral closure of $A$ in $A_K$; in view of
claim \ref{cl_first-after-K} we are reduced to showing that
$A=A^\nu$. By proposition \ref{prop_integr-closure-grad}(ii),
$A^\nu$ is $\Gamma$-graded.
Suppose now that $x\in A^\nu$; we need to show that $x\in A$,
and we may assume that $x\in\gr_\gamma A^\nu$ for some
$\gamma\in\Gamma$. Hence, let
$$
x^n+a_1x^{n-1}+\cdots+a_n=0
\qquad
\text{with $a_1,\dots,a_n\in A$}
$$
be an integral equation for $x$ over $A$. If we replace each
$a_i$ by its homogeneous component in degree $i\gamma$, we
still obtain an integral equation for $x$, so we may assume
that $a_i\in\gr_{i\gamma}A$ for every $i=1,\dots,n$.
Then, (MA2) implies that, for every $i=1,\dots,n$ there exist
$u_i\in K^+$ and $b_i\in\gr_\gamma A$, such that $|u_i|=1$ and
$a_i=u_ib_i^i$. Also, from (MA1) we deduce that there exists
$b\in\gr_\gamma A$ such that $b_1,\dots,b_n\in K^+b$.
Set $y:=b^{-1}x\in\gr_0A\otimes_{K^+}K$; clearly $y$ is integral
over the subring $\gr A_0$. Lastly, (MA3) shows that $y\in K^+$,
whence $x\in\gr_\gamma A$, which proves the contention.
\end{proof}

\begin{proposition}\label{prop_need-it-now}
Let $(A,\Gamma)$ be a model $K^+$-algebra, and suppose
that $\Gamma_{\!\!K}$ is divisible. Then $(A,\Gamma)$ is
the filtered union of its small model $K^+$-subalgebras.
\end{proposition}
\begin{proof} In view of lemma \ref{lem_normal-norman}(i), we
see that $(A,\Gamma)$ is the colimit of the filtered system
of its subobjects $(A_{\Delta[1/p]},\Delta[1/p])$, for $\Delta$
ranging over the fine and saturated submonoids of $\Gamma$.
We may then assume from start that $\Gamma=\Gamma_{\!0}[1/p]$
for some fine and saturated submonoid $\Gamma_{\!0}$.

Next, in view of remark \ref{rem_decompose-in-B}(i), we may consider
separately the cases where $\Gamma_{\!0}^\gp$ is a torsion-free
abelian group, and where $\Gamma=\Gamma_{\!0}$ is a finite abelian
group.

Suppose first that $\Gamma_{\!0}^\gp$ is torsion-free, and let
$\underline\gamma:=(\gamma_1,\dots,\gamma_n)$ be a finite system
of generators for $\Gamma_{\!0}$. For every $i=1,\dots,n$, choose
$a_i\in A^*_{\gamma_i}$, and let
$B(\underline\gamma,\underline a)\subset A$ be the $K^+$-subalgebra
generated by $\underline a:=(a_1,\dots,a_n)$.
Clearly the grading of $A$ induces a $\Gamma_{\!0}$-grading on
$B(\underline\gamma,\underline a)$, and
$\gr_\beta B(\underline\gamma,\underline a)$ is a finitely
generated $K^+$-module for every $\beta\in\Gamma_{\!0}$
(proposition \ref{prop_four-year-later}(ii)); by virtue of
(MA1), we know that $\gr_\beta B(\underline\gamma,\underline a)$
is then even a free rank one $K^+$-module, for every
$\beta\in\Gamma_{\!0}$. Furthermore, $\underline a$ generates a
fine submonoid of $A^*$, so -- by proposition \ref{prop_flattening}(ii) --
there exists an integer $k>0$ such that
$$
(\gr_{k\beta} B(\underline\gamma,\underline a))^n=
\gr_{nk\beta}B(\underline\gamma,\underline a)
\qquad
\text{for every $\beta\in\Gamma_{\!0}$ and every integer $n>0$}.
$$
Now, for any $\beta\in\Gamma$, let $t>0$ be an integer such
that $p^t\beta\in\Gamma_0$, and pick a generator $b$ of
the $K^+$-module $\gr_{kp^t\beta}B(\underline\gamma,\underline a)$;
in light of (MA2) we may find $x\in K^+$ and $c\in\gr_\beta A$,
such that $b=c^{kp^t}x$. Since $\Gamma_{\!\!K}$ is divisible, we
may write $x=y^{kp^t}u$ for some $y,u\in K^+$ with $|u|=1$.

It is easily seen that the $K^+$-submodule of $\gr_\beta A$
generated by $cy$ does not depend on the choices of $c$, $y$,
$t$ and $k$; hence we denote
$\gr_\beta C(\underline\gamma,\underline a)$ this submodule;
a simple inspection shows that the resulting $\Gamma$-graded
$K^+$-module
$$
C(\underline\gamma,\underline a):=
\bigoplus_{\beta\in\Gamma}\gr_\beta C(\underline\gamma,\underline a)
$$
is actually a $\Gamma$-graded $K^+$-subalgebra of $A$,
for which (MA2) holds, and therefore
$(C(\underline\gamma,\underline a),\Gamma)$ is a small model
algebra. Lastly, it is clear that the family of all such
$(C(\underline\gamma,\underline a),\Gamma)$, for $\underline\gamma$
ranging over all finite sets of generators of $\Gamma$, and
$\underline a$ ranging over all the finite sequences of elements
of $A$ as above, form a cofiltered system of subobjects, whose
colimit is $(A,\Gamma)$. This concludes the proof of the proposition
in this case.

Next, suppose that $\Gamma=\Gamma_{\!0}$ is a finite abelian
group. In view of remark \ref{rem_tensor-prod-in-A}(iii), we
are reduced to the case where $\Gamma$ is cyclic, say
$\Gamma=\Z/n\Z$. Fix any generator $\gamma$ of $\Gamma$,
and for every $a\in\gr_\gamma A$, let $C(a)$ be the $K^+$-subalgebra
of $A$ generated by $a$; it is easily seen that the colimit
of the filtered family $((C(a),\Gamma)~|~a\in\gr_\gamma A)$
equals $(A,\Gamma)$. (Details left to the reader.)
\end{proof}

\sset\subsubsection{}\label{subsec_chart-for purity}
Suppose now that $K$ is algebraically closed, let $(B,\Delta)$
be a small model $K^+$-algebra, with $\Delta^\gp$ torsion-free,
and write $\Delta=\Delta_0[1/p]$ for some fine and saturated
submonoid $\Delta_0$ such that the $K^+$-algebra $B_0:=B_{\Delta_0}$
is finitely generated. Let us set
$$
\Delta_n:=\{\gamma\in\Delta~|~p^n\gamma\in\Delta_0\}
\qquad
B_n:=B_{\Delta_n}
\qquad
Y_n:=\Spec\,B_n
\qquad
\text{for every $n\in\N$}.
$$
We wish to construct a ladder of log schemes
\set\begin{equation}\label{eq_log-ladder}
{\diagram \cdots \ar[r]^-{g_{n+1}} &
(Y_{n+1},\underline M{}_{n+1}) \ar[r]^-{g_n} \ar[d]_{\phi_{n+1}} &
(Y_n,\underline M{}_n) \ar[r]^-{g_{n-1}} \ar[d]_{\phi_n} &
\cdots \ar[r]^-{g_0} & (Y_0,\underline M{}_0) \ar[d]^{\phi_0} \\
\cdots \ar[r]^-{h_{n+1}} & (S,\underline N{}_{n+1}) \ar[r]^-{h_n} &
(S,\underline N{}_n) \ar[r]^-{h_{n-1}} &
\cdots \ar[r]^-{h_0} & (S,\underline N{}_0)
\enddiagram}
\end{equation}
such that, for every $n\in\N$ :
\begin{itemize}
\item
$\phi_n$ is smooth and saturated, and $\underline M{}_n$ admits
a chart, given by a $\Delta_n$-graded fine monoid $P^{(n)}$,
such that $P^{(n)}_0$ is sharp, and the inclusion map
$P^{(n)}_0\to P^{(n)}$ is flat and saturated, and gives
a chart for $\phi_n$
\item
the morphism of schemes underlying $h_n$ (resp. $g_n$) is the
identity of $S$ (resp. is induced by the inclusion
$B_n\subset B_{n+1}$)
\item
the morphism $g_n$ admits a chart, given by an injective
map $P^{(n)}\to\Delta_n\times_{\Delta_{n+1}}P^{(n+1)}$
of $\Delta_n$-graded monoids, whose restriction
$P^{(n)}_0\to P^{(n+1)}_0$ gives a chart for $h_n$.
\end{itemize}
This will be achieved in several steps, as follows :

$\bullet$\ \
First, let $b_1,\dots,b_k$ be a finite system of generators
for $B_0$. By remark \ref{rem_good-luck-wt}(ii), we have a
decomposition $B^*_0=(K^+)^\times\oplus C$ for some submonoid
$C\subset B^*_0$, and we may suppose that each $b_i$ lies in
$C$. Let $P^{(0)}\subset B^*_0$ be the submonoid generated
by $b_1,\dots,b_k$. The restriction of the surjection
$B^*_0\to\Delta_0$ is then still a surjection
$\pi:P^{(0)}\to\Delta_0$. Notice that, for every
$x\in\Ker\,\pi^\gp$, we have either $x\in K^+$ or
$x^{-1}\in K^+$; especially, we may find a finite system
$x_1,\dots,x_n$ of generators for $\Ker\,\pi^\gp$, such
that $x_i\in K^+$ for every $i=1,\dots,n$. Now, let
$\Sigma\subset P^{(0)\gp}$ be the submonoid generated by
$x_1,\dots,x_n$; after replacing $P^{(0)}$ by
$P^{(0)}\cdot\Sigma$, we may assume that
$\Ker\,\pi^\gp=P_0^{(0)\gp}$, where
$P^{(0)}_0:=\Ker\,\pi^\gp\cap P^{(0)}\subset K^+$ is a fine
submonoid.

$\bullet$\ \
Next, by theorem \ref{th_flattening-th}, we may find a finitely
generated submonoid $\Sigma'\subset P_0^{(0)\gp}\cap K^+$ such that
the induced morphism $P^{(0)}_0\cdot\Sigma'\to P^{(0)}\cdot\Sigma'$
is flat. Clearly
$$
P^{(0)}_0\cdot\Sigma'=\Ker\,\pi^\gp\cap(P^{(0)}\cdot\Sigma')
$$
hence, we may replace $P^{(0)}$ by $P^{(0)}\cdot\Sigma'$, and
assume that the morphism $P^{(0)}_0\to P^{(0)}$ is also flat.

$\bullet$\ \
We claim that the map $P^{(0)}_0\to P^{(0)}$ is also saturated.
We shall apply the criterion of proposition
\ref{prop_crit-saturated} : indeed, for given $\gamma\in\Delta_0$
and integer $n>0$, let $x$ (resp. $y$) be a generator of the
$P^{(0)}_0$-module $P^{(0)}_\gamma$ (resp. $P^{(0)}_{n\gamma}$).
Then $x$ (resp. $y$) is also a generator of the $K^+$-module
$\gr_\gamma B$ (resp. $\gr_{n\gamma}B$), and (MA2) implies that
there exists $u\in(K^+)^\times$ such that $uy=x^n$; but since
$x,y\in C$, we must have $u=1$ therefore
$(P^{(0)}_\gamma)^n=P^{(0)}_{n\gamma}$, as required.

$\bullet$\ \
Next, we claim that the induced map of $\Delta_0$-graded
$K^+$-algebras
\set\begin{equation}\label{eq_so-simple}
P^{(0)}\otimes_{P^{(0)}_0}K^+\to B_0
\end{equation}
is an isomorphism. Indeed, since $P^{(0)}$ contains a set of
generators for the $K^+$-algebra $B$, the map \eqref{eq_so-simple}
is obviously surjective. However, for every $\gamma\in\Delta_0$,
the $P^{(0)}_0$-module $P^{(0)}_\gamma$ is free of rank one
(remark \ref{rem_why-not}(iv)),
therefore $(P^{(0)}\otimes_{P^{(0)}_0}K^+)_\gamma$ is a free
$K^+$-module of rank one. It follows easily that
\eqref{eq_so-simple} is also injective.

$\bullet$\ \
Now, denote by $\underline N{}_0$ (resp. $\underline M{}_0$)
the fine log structure on $S$ (resp. on $X_0$) deduced from the
inclusion map $P^{(0)}_0\to K^+$ (resp. $P^{(0)}\to B_0$).
By lemma \ref{lem_integr-flat}(iv), the inclusion
$P^{(0)}_0\to P^{(0)}$ yields a chart for a morphism
$\phi_0:(Y_0,\underline M{}_0)\to(S,\underline N{}_0)$
that is saturated, as sought. Lastly, since \eqref{eq_so-simple}
is an isomorphism, theorem \ref{th_charact-smoothness} shows
that $\phi_0$ is also smooth.

$\bullet$\ \
Next, according to remark \ref{rem_good-luck-wt}(iii), we have
a compatible system of decompositions
$$
B_n^*=(K^+)^\times\oplus C^{(n)}
\qquad
\text{for every $n\in\N$}
$$
with $C^{(0)}=C$, and such that the $p$-Frobenius induces an
isomorphism $\tau_n:C^{(n+1)}\isom C^{(n)}$ for every $n\in\N$.
Hence, define inductively an increasing sequence of submonoids
$$
P^{(0)}\subset P^{(1)}\subset P^{(2)}\subset\cdots\subset B^*
\quad\text{by letting}\quad
\text{$P^{(n+1)}:=\tau_n^{-1}P^{(n)}$ for every $n\in\N$}.
$$
Clearly, the grading of $B^*$ restricts to a $\Delta_n$-grading
on $P^{(n)}$, and induces a $\Delta$-grading on $P$. We deduce
isomorphisms of $\Delta_{n+1}$-graded monoids :
\set\begin{equation}\label{eq_change-grading}
P^{(n+1)}\isom\Delta_{n+1}\times_{\Delta_n}P^{(n)}
\qquad
\text{for every $n\in\N$.}
\end{equation}
Especially, the induced maps of monoids $P^{(n)}_0\to P^{(n)}$
are still flat and saturated, and $P^{(n)}_0$ is still sharp
(remark \ref{rem_unique-generator}). The inclusion maps
$P^{(n)}\to B_n$ and $P^{(n)}_0\to K^+$ determine an
isomorphism of $\Delta_n$-graded $K^+$-algebras
\set\begin{equation}\label{eq_simple-enough}
P^{(n)}\otimes_{P^{(n)}_0}K^+\isom B_n
\qquad
\text{for every $n\in\N$}
\end{equation}
as well as fine log structures $\underline M{}_n$ on $Y_n$,
and $\underline N{}_n$ on $S$, whence a morphism of log schemes
$\phi_n:(Y_n,\underline M{}_n)\to(S,\underline N{}_n)$ which is
again smooth and saturated. This completes the construction of
\eqref{eq_log-ladder}.

\begin{remark}\label{rem_unique-generator}
(i)\ \
With the notation of \eqref{subsec_chart-for purity}, notice that,
by construction, $P^{(n)\gp}\subset C^{(n)\gp}$; especially, since
$\Delta^\gp$ is torsion-free, the same holds for $P^{(n)\gp}$.
Likewise, the construction shows that the induced map
$P^{(n)\gp}_0\to\Gamma_{\!\!K}$ is always injective; especially,
$P^{(n)\gp}_0$ is a free abelian group of finite rank, and
$P^{(n)}_0$ is sharp. It follows that the $P^{(n)}_0$-module
$P^{(n)}_\gamma$ admits a unique generator $g_\gamma$, for every
$\gamma\in\Delta_n$, and the saturation condition for the inclusion
$P^{(n)}_0\to P^{(n)}$ translates as the system of identities :
$$
g_\gamma^k=g_{k\gamma}
\qquad
\text{for every $k\in\N$ and every $\gamma\in\Delta_n$}.
$$

(ii)\ \
Moreover, set $P:=\bigcup_{n\in\N}P^{(n)}$. Notice that
the colimit of the maps \eqref{eq_simple-enough} is an
isomorphism of $\Delta$-graded $K^+$-algebras
\set\begin{equation}\label{eq_simple-colim-enough}
P\otimes_{P_0}K^+\isom B.
\end{equation}
Furthermore, since the natural map $P^{(n)}/P^{(n)}_0\to\Delta$
is an isomorphism for every $n\in\N$, we deduce that the
grading of $P$ induces an isomorphism :
\set\begin{equation}\label{eq_state-this-iso}
P/P_0\isom\Delta.
\end{equation}

(iii)\ \
We also obtain a commutative diagram of log schemes :
\set\begin{equation}\label{eq_more-esses}
{\diagram
\SS(B_n,\Delta_n) \ar[r] \ar[d] &
(Y_n,\underline M{}_n) \ar[d]^{\phi_n} \\
\SS(K^+) \ar[r] & (S,\underline N{}_n)
\enddiagram}
\qquad
\text{for every $n\in\N$}
\end{equation}
whose left vertical arrow is the saturated morphism
\eqref{eq_ss-aturated}, and whose bottom (resp. top) arrow
is the identity of $S$ (resp. of $Y_n$) on the underlying
schemes, and is induced by the inclusion map $P^{(0)}_n\to K^*$
(resp. $P^{(n)}\to B^*_n$). It is easily seen that
\eqref{eq_more-esses} is cartesian : indeed, since the
functor \eqref{eq_another-left-adj} is right exact, it
suffices to check that the natural map
$$
P^{(n)}\otimes_{P^{(n)}_0}K^*\to B^*_n
$$
is an isomorphism, which is clear from \eqref{eq_simple-enough}.

(iv)\ \
Notice that the inclusion map $i_0:P^{(0)}_0\to P^{(0)}$
is also local; then corollary \ref{cor_general-deco} yields a
decomposition
$$
\theta_0:P^{(0)}\isom P^{(0)\times}\oplus P^{(0)\sharp}
$$
such that $\theta_0\circ i_0=\iota_0\circ\lambda^\sharp_0$,
where $\iota_0:P^{(0)\sharp}\to P^{(0)\times}\oplus P^{(0)\sharp}$
is the natural inclusion map. By means of the isomorphisms
\eqref{eq_change-grading}, we may then inductively construct
decompositions
$$
\theta_n:P^{(n)}\isom P^{(n)\times}\oplus P^{(n)\sharp}
\qquad
\text{for every $n\in\N$}
$$
fitting into a commutative diagram
$$
\xymatrix{
P^{(n)}_0 \ar[rr]^-{i_n} \ar[d]_{i_n^\sharp} & &
P^{(n)} \ar[rr]^-{j_n} \ar[d]_{\theta_n} & &
P^{(n+1)} \ar[d]^{\theta_{n+1}} \\
P^{(n)\sharp} \ar[rr]^-{\iota_n} & &
P^{(n)\times}\oplus P^{(n)\sharp}
\ar[rr]^-{j_n^\times\oplus j_n^\sharp} & &
P^{(n+1)\times}\oplus P^{(n+1)\sharp}
}$$
where $i_n$, $j_n$ and $\iota_n$ are the natural inclusion
maps. Now, set
$$
B'_n:=K^+[P^{(n)\times}]
\qquad
B''_n:=P^{(n)\sharp}\otimes_{P^{(n)}_0}K^+
\qquad
\text{for every $n\in\N$}
$$
and let $\Delta'_n$ (resp. $\Delta''_n$) be the image of
$P^{(n)\times}$ (resp. of $\theta^{-1}_nP^{(n)\sharp}$)
in $\Delta_n$; since the grading of $P^{(n)}$ induces an
isomorphism $P^{(n)}/P^{(n)}_0\isom\Delta_n$, we see that
$\Delta_n=\Delta'_n\oplus\Delta''_n$, and by construction,
for every $n\in\N$ we have isomorphisms
$$
(B_n,\Delta_n)\isom(B'_n,\Delta'_n)\otimes(B''_n,\Delta''_n)
$$
of model $K^+$-algebras, which identify the inclusions
$B_n\to B_{n+1}$ with the tensor product of the induced
inclusion maps $B'_n\to B'_{n+1}$ and $B''_n\to B''_{n+1}$.
Furthermore, set $Y'_n:=\Spec\,B'_n$ and $Y''_n:=\Spec\,B''_n$
for every $n\in\N$. The induced map of monoids
$P^{(n)\sharp}\to B''_n$ determines a fine log structure
$\underline M{}''_n$ on $Y''_n$, the map $i^\sharp_n$ gives
a chart for a smooth and saturated morphism of log schemes
$$
\phi''_n:(Y''_n,\underline M{}''_n)\to(S,\underline N{}_n)
$$
and there follows an isomorphism of
$(S,\underline N{}_n)$-schemes (lemma \ref{lem_check-iso})
\set\begin{equation}\label{eq_decompose-ladder}
(Y_n,\underline M{}_n)\isom
Y'_n\times_S(Y''_n,\underline M{}''_n)
\qquad
\text{for every $n\in\N$}.
\end{equation}

(v)\ \
Starting with \eqref{subsec_extend-to-Q}, we shall consider
the strict henselization $B^\sh$ of $B$ at a given geometric
point $\bar x$ of $\Spec\,B\otimes_{K^+}\kappa$. In this
situation, let $x\in Y:=\Spec\,B$ be the support of $\bar x$,
and $\fp_x\subset B$ the corresponding prime ideal. Let also
$P$ be as in (ii), denote by $\beta:P\to B$ the natural
map deduced from \eqref{eq_simple-colim-enough}, and set
$\fp:=\beta^{-1}\fp_x$.
Moreover, set $\fp_n:=\fp\cap P^{(n)}$ and
$Q^{(n)}:=P^{(n)}_{\fp_n}$ for every $n\in\N$;
clearly $\fp_n\subset\fp_{n+1}$ (so the isomorphism
\eqref{eq_change-grading} maps $\fp_{n+1}$ onto $\fp_n$),
and therefore
$$
Q^{(n)}\subset Q^{(n+1)}
\qquad
\text{for every $n\in\N$}.
$$
Furthermore, the image of $Q^{(n)}$ in $\Delta^\gp_n$
is a localization $\Gamma_n$ of $\Delta_n$, especially it
is still saturated, and the maps \eqref{eq_change-grading}
extend to isomorphisms of $\Gamma_{n+1}$-graded monoids
$$
Q^{(n+1)}\isom\Gamma_{n+1}\times_{\Gamma_n}Q^{(n)}
\qquad
\text{for every $n\in\N$}.
$$
It is also clear that the $p$-Frobenius of $\Gamma_{n+1}$
factors through an isomorphism $\Gamma_{n+1}\isom\Gamma_n$,
for every $n\in\N$. Set $F_n:=P^{(n)}\setminus\fp_n$ for
every $n\in\N$; notice that, by construction,
$F_n\cap P^{(n)}_0=\{1\}$, hence $F_n\cap P^{(n)}_\gamma$
is either empty or else it contains exactly one element,
namely the generator $g_\gamma$ of $P^{(n)}_\gamma$, by
virtue of (i).
It follows that $Q^{(n)}_0=P^{(n)}_0$ for every $n\in\N$,
and therefore the inclusion map $Q^{(n)}_0\to Q^{(n)}$
is still flat and saturated (lemma \ref{lem_little}(ii)).
Let $\Gamma:=\bigcup_{n\in\N}\Gamma_n$; we conclude that
$Q:=P_\fp$ is a $\Gamma$-graded monoid with $Q_0=P_0$, and
we may define
$$
B_\fp:=Q\otimes_{Q_0}K^+.
$$
The foregoing shows that $(B_\fp,\Gamma)$ is still a small
model algebra, and we also obtain a ladder of log schemes
with the properties listed in \eqref{subsec_chart-for purity} :
namely, set
$B_{\fp,n}:=B_{\fp,\Gamma_n}=Q^{(n)}\otimes_{Q^{(n)}_0}K^+$
for every $n\in\N$, and endow $Y_{\fp,n}:=\Spec\,B_{\fp,n}$
with the log structure $\underline M{}_{\fp,n}$ determined
by the induced map $Q^{(n)}\to B_{\fp,n}$. Clearly the
geometric point $\bar x$ lifts uniquely to a geometric
point $\bar x_\fp$ of $Y_\fp:=\Spec\,B_\fp$, and the
localization map $B\to B_\fp$ induces an isomorphism
$$
Y_\fp(\bar x_\fp)\isom Y(\bar x).
$$
Hence, for the study of the scheme $Y(\bar x)$, it shall
be usually possible to replace the original small algebra
$B$ by its localization $B_\fp$ thus constructed. In so
doing, we gain one more property: indeed, notice that
the new chart $Q^{(n)}\to B_{\fp,n}$ is local at the
geometric point $\bar x$, for every $n\in\N$. Now, choose
a compatible system of decompositions
$$
(Y_{\fp,n},\underline M{}_{\fp,n})\isom
Y'_{\fp,n}\times_S(Y''_{\fp,n},\underline M{}''_{\fp,n})
$$
as in (iv), and denote by $\bar x_{\fp,n}$ (resp.
$\bar x{}''_{\fp,n}$) the image of $\bar x_\fp$ in
$Y_{\fp,n}$ (resp. in $Y''_{\fp,n}$), for every $n\in\N$. By
inspecting the construction, it is easily seen that the
chart $Q^{(n)\sharp}\to B''_{\fp,n}:=
Q^{(n)\sharp}\otimes_{Q_0^{(n)}}K^+$ is also local at
$\bar x{}''_{\fp,n}$.
\end{remark}

\sset\subsubsection{}\label{subsec_extend-to-Q}
Let us return to the situation of \eqref{subsec_chart-for purity}.
Fix a geometric point $\bar x$ of $Y:=\Spec\,B$, localized
on $Y\times_S\Spec\,\kappa$, and for every $n\in\N$, let
$\bar x_n$ be the image of $\bar x$ in $Y_n$. Set
$$
B^\sh:=\cO_{Y(\bar x),\bar x}
\qquad\text{and}\qquad
B_n^\sh:=\cO_{Y_n(\bar x_n),\bar x_n}
\quad
\text{for every $n\in\N$}.
$$
For the next step, we shall apply the method of
\eqref{subsec_last-step}, to construct a normalized length
function for $B^\sh$-modules. To this aim, we need an
auxiliary model $K^+$-algebra, defined as follows.
First, notice that the grading $P\to\Delta$ extends to
a morphism of monoids
$$
\pi_\Q:P_\Q\to\Delta_\Q
$$
(notation of \eqref{subsec_from-con-to-mon}); {\em i.e.}
$P_\Q$ is a $\Delta_\Q$-graded monoid. Since $K^\times$ is
a divisible group, the inclusion map $P_0\to K^*$ extends
to a group homomorphism $P^\gp_{\Q,0}\to K^\times$, and it
is easily seen that the latter restricts to a morphism of monoids
$$
P_{\Q,0}\to K^*.
$$
So finally, we may set
$$
A:=P_\Q\otimes_{P_{\Q,0}}K^+.
$$
Taking into account \eqref{eq_simple-colim-enough}, we deduce
an isomorphism of $\Delta$-graded $K^+$-algebras 
\set\begin{equation}\label{eq_deduce-this}
B\isom A_\Delta.
\end{equation}
Moreover, arguing as in remark \ref{rem_good-luck-wt}(iii),
we get an isomorphism
\set\begin{equation}\label{eq_same-as-ever}
K[\Delta_\Q]\isom A_K
\end{equation}
fitting into a commutative diagram of $K$-algebras
$$
\xymatrix{ K[\Delta] \ar[r]^-\sim \ar[d] & B_K \ar[d] \\
K[\Delta_\Q] \ar[r]^-\sim & A_K
}$$
whose top horizontal arrow is the isomorphism of remark
\ref{rem_decompose-in-B}(v), and whose right (resp. left)
vertical arrow is deduced from \eqref{eq_deduce-this}
(resp. is induced by the natural inclusion
$\Delta\subset\Delta_\Q$).

\begin{lemma}\label{lem_extend-to-Q}
With the notation of \eqref{subsec_extend-to-Q}, we have :
\begin{enumerate}
\item
$P_{\Q,\gamma}$ is a free $P_{\Q,0}$-module
of rank one, for every $\gamma\in\Delta_\Q$.
\item
$P_{\Q,0}$ is a sharp monoid, and the inclusion map
$P_{\Q,0}\to P_\Q$ is flat and saturated.
\end{enumerate}
\end{lemma}
\begin{proof}(i): We can write $P_\Q$ as the colimit of
the system of monoids
$$
(P_{[n]};\ \mu_{P,n,m}:P_{[n]}\to P_{[nm]}~|~n,m\in\N)
$$
such that $P_{[n]}:=P^{(0)}$ for every $n\in\N$, and $\mu_{P,n,m}$
is the $m$-Frobenius map of $P^{(0)}$, for every $n,m\in\N$. Likewise,
we may write $\Delta_\Q$ as colimit of a system of Frobenius
endomorphisms $(\Delta_{0,[n]};\ \mu_{\Delta,n,m}~|~n,m\in\N)$,
and $\pi_\Q$ is the colimit of the corresponding system of maps
$(\pi_{[n]}:P_{[n]}\to\Delta_{0,[n]}~|~n\in\N)$ where
$\pi_{[n]}:=\pi$ for every $n\in\N$.
Hence, for any $\gamma\in\Delta_0$ there exists $n\in\N$ such
that $\gamma$ is the image of some $\gamma_n\in\Delta_{0,[n]}$,
and
$$
P_{\Q,\gamma}=\colim_{k\in\N}P_{[nk],k\gamma_n}=
\colim_{k\in\N}P_{k\gamma_n}.
$$
However, we know that $P_{k\gamma_n}$ is free of rank one,
generated by a unique element $g_{k\gamma_n}$; moreover,
the transition maps $\mu_{P,n,k}$ send $g_{\gamma_n}$
onto $g_{k\gamma_n}$, for every $k\in\N$ (remark
\ref{rem_unique-generator}). This implies that $P_{\Q,\gamma}$
is generated by the image of $g_{\gamma_n}$, whence (i).

(ii): The flatness follows from (i). Since the inclusion
map $P^{(0)}_0\to P^{(0)}$ is saturated, it is easily seen
that the same holds for the inclusion $P_{\Q,0}\to P_\Q$.
Lastly, the sharpness of $P_{\Q,0}$ is likewise deduced
from the sharpness of $P^{(0)}_0$ (details left to the reader).
\end{proof}

\sset\subsubsection{}\label{subsec_define-f_A}
Lemma \ref{lem_extend-to-Q}(ii) implies that the pair
$(A,\Delta_\Q)$ fulfills axiom (MA2), and we have thus our
sought auxiliary model $K^+$-algebra.
Since $P_{\Q,0}$ is sharp, the $P_{\Q,0}$-module $P_{\Q,\gamma}$
admits a unique generator $g_\gamma$, for every $\gamma\in\Delta_\Q$.
The image $g_\gamma\otimes 1$ of $g_\gamma$ in $A$ is a generator
of the direct summand $A_\gamma$, which is a free $K^+$-module of
rank one; hence, there exists a unique $a_\gamma\in K$ such that
$\gamma\otimes a_\gamma$ gets mapped to $g_\gamma\otimes 1$, under
the isomorphism \eqref{eq_same-as-ever}.

After choosing an order-preserving isomorphism
\eqref{eq_law-and-order}, we may define a function
$$
f_A:\Delta_\Q\to\R
\quad : \quad
\gamma\mapsto\log|a_\gamma|.
$$
The inclusion $A_\gamma\cdot A_\delta\subset A_{\gamma+\delta}$
translates as the inequality
$$
f_A(\gamma)+f_A(\delta)\geq f_A(\gamma+\delta)
\qquad
\text{for every $\gamma,\delta\in\Delta_\Q$}.
$$
Likewise, the saturation condition of axiom (MA2) translates
as the identity
$$
f_A(n\gamma)=n\cdot f_A(\gamma)
\qquad
\text{for every $\gamma\in\Delta_\Q$ and every $n\in\N$}.
$$
We also fix a (Banach) norm
$$
\Delta^\gp_\R\to\R
\qquad
\gamma\mapsto\Vert\gamma\Vert.
$$
We can then state the following

\begin{lemma}\label{lem_liner-subdivide}
With the notation of \eqref{subsec_define-f_A}, we have :
\begin{enumerate}
\item
There exists a $\Delta_0^\gp$-rational subdivision $\Theta$
of the convex polyhedral cone $(\Delta^\gp_\R,\Delta_\R)$,
such that
$$
f_A(\gamma+\delta)=f_A(\gamma)+f_A(\delta)
\qquad
\text{for every $\sigma\in\Theta$ and every
$\gamma,\delta\in\sigma\cap\Delta_\Q$}.
$$
\item
Especially, the function $f_A$ is of\/ {\em Lipschitz type},
{\em i.e.} there exists a real constant $C_A>0$ such that
$$
|f_A(\gamma)-f_A(\gamma')|\leq C_A\cdot\Vert\gamma-\gamma'\Vert
\qquad
\text{for every $\gamma,\gamma'\in\Delta^\gp_\Q$}.
$$
\end{enumerate}
\end{lemma}
\begin{proof} (Here, $\Delta^\gp_\R:=\Delta^\gp_0\otimes_\Z\R$
and $\Delta_\R\subset\Delta^\gp_\R$ is the cone spanned by
$\Delta_0$.) Combining proposition \ref{prop_divide-et-imp}
and lemma \ref{lem_not-so-elegant}, we find a
$\Delta_0^\gp$-rational subdivision $\Theta$ of $\Delta_\R$
such that
$$
P_{\Q,\gamma+\delta}=P_{\Q,\gamma}+P_{\Q,\delta}
\qquad
\text{for every $\sigma\in\Theta$ and every
$\gamma,\delta\in\sigma\cap\Delta_\Q$}.
$$
This means that $g_\gamma\cdot g_\delta=g_{\gamma+\delta}$
for every $\sigma\in\Theta$ and every
$\gamma,\delta\in\sigma\cap\Delta_\Q$, whence (i).

(ii): Clearly, for every $\sigma\in\Theta$ we may find
a constant $C_\sigma$ such that the stated inequality
holds -- with $C_A$ replaced by $C_\sigma$ -- for every
$\gamma,\gamma'\in\sigma\cap\Delta_\Q$. It is easily
seen that $C_A:=\max(C_\sigma~|~\sigma\in\Theta)$ will do
(details left to the reader).
\end{proof}

\sset\subsubsection{}\label{subsec_tenero-giaco}
Fix $n\in\N$, pick a subdivision $\Theta$ of
$(\Delta^\gp_\R,\Delta_\R)$ as in lemma \ref{lem_liner-subdivide},
and let $\Theta^s\subset\Theta$ be the subset of all
$\sigma\in\Theta$ that span $\Delta^\gp_\R$. Also, for
every $\gamma\in\Delta^\gp_\Q$, denote by $[\gamma]$
the class of $\gamma$ in $\Delta^\gp_\Q/\Delta^\gp_n$.
For any subset $\Sigma\subset\Delta_\R^\gp$ we let
$$
A_\Sigma:=A_{\Sigma\cap\Delta_\Q}
$$
(where the right-hand side is defined as in
\eqref{subsec_not-so-bad}). In light of lemma
\ref{lem_lem_top-dim-subdi} we obtain a decomposition
of $B_n=A_{\Delta_n}=A_{[0]}$ as sum of $K^+$-subalgebras :
$$
B_n=\sum_{\sigma\in\Theta^s}A_{[0]\cap\sigma}
$$
and a corresponding decomposition of the $B_n$-module
$A_{[\gamma]}$ as sum of $A_{[0]\cap\sigma}$-modules
$$
A_{[\gamma]}=
\sum_{\sigma\in\Theta^s}A_{[\gamma]\cap\sigma}
\qquad
\text{for every $\gamma\in\Delta^\gp_\Q$}.
$$

\begin{lemma}\label{lem_choose-int-N}
With the notation of \eqref{subsec_tenero-giaco} and
\eqref{subsec_chart-for purity}, there exists $N\in\N$ such
that for every $n\in\N$ and every $\gamma\in\Delta^\gp_\Q$,
the $B_n$-module $A_{[\gamma]}$ admits a system of generators
of cardinality at most $N$.
\end{lemma}
\begin{proof} According to proposition \ref{prop_Gordon}(ii),
for every $\sigma\in\Theta^s$ and every $\gamma\in\Delta^\gp_\Q$,
the subset $\Sigma_{\sigma,\gamma}:=\Delta^\gp_n\cap(\sigma-\gamma)$
is a finitely generated $(\Delta_n^\gp\cap\sigma)$-module; hence,
let us fix a finite system of generators $G_{\sigma,\gamma}$
for $\Sigma_{\sigma,\gamma}$. In light of lemma
\ref{lem_liner-subdivide}(i) we see that the finite set
$$
G'_{\sigma,\gamma}:=
\{g_{\gamma+\delta}\otimes 1~|~\delta\in G_{\sigma,\gamma}\}
\subset A
$$
generates the $A_{[0]\cap\sigma}$-module $A_{[\gamma]\cap\sigma}$.
Consequently, the finite set
$\bigcup_{\sigma\in\Theta^s}G'_{\sigma,\gamma}$
generates the $B_n$-module $A_{[\gamma]}$.
On the other hand, since $\cS_{\!\Delta^\gp_n,\sigma}$ is
a finitely generated $\Delta^\gp_n$-module (proposition
\ref{prop_linear-constr-part}(ii))
it is clear that the cardinality of $G_{\sigma,\gamma}$ is
bounded by a constant $N_\sigma$ that is independent of
$\gamma$, and then $A_{[\gamma]}$ is generated by at most
$N_n:=\sum_{\sigma\in\Theta^s}N_\sigma$ elements.
It remains to show that the estimate for $N_n$ is independent
of $n$. However, notice that, for every $\sigma\in\Theta^s$,
the automorphism of $\Delta^\gp_\Q$ given by multiplication
by $p^n$, induces a natural bijection
$$
\cS_{\!\Delta^\gp_n,\sigma}\isom\cS_{\!\Delta^\gp_0,\sigma}
$$
that sends each $\Delta_n$-module $\Delta^\gp_n\cap(\sigma-v)$
onto the $\Delta_0$-module
$$
\Delta^\gp_0\cap(\sigma-p^nv)=
\Delta_0\otimes_{\Delta_n}(\Delta^\gp_n\cap(\sigma-v))
$$
(where the extension of scalars $\Delta_n\to\Delta_0$ is
the isomorphism given by the rule : $\gamma\mapsto p^n\gamma$
for every $\gamma\in\Delta_n\cap\sigma$). Thus, we see that
$N:=N_0$ will already do.
\end{proof}

\sset\subsubsection{}
Let us choose $N$ as in lemma \ref{lem_choose-int-N}, and
endow the set $\cM_N(B^a_n)$ with the uniform structure
defined in \eqref{subsec_uniform} (relative to the standard
setup attached to $K^+$); we consider the mapping
$$
\Delta^\gp_\Q\to\cM_N(B^a_n)
\qquad
\gamma\mapsto A_{[\gamma]}^a.
$$

\begin{lemma}\label{lem_lot-of-notate}
Keep the notation of \eqref{subsec_for-almost-lengths} and
\eqref{subsec_uniform}, and let also $C_A>0$ be the constant
appearing in lemma {\em\ref{lem_liner-subdivide}(ii)}. Then :
$$
(A_{[\gamma]}^a,A_{[\gamma']}^a)
\in E_{C_A\cdot\Vert\gamma-\gamma'\Vert\cdot 2}
$$
for every $\Sigma\in\cS_{\!\Delta^\gp_n,\Delta_\R}$ and every
$\gamma,\gamma'\in\Omega(\Delta_\R,\Sigma)\cap\Delta^\gp_\Q$.
\end{lemma}
\begin{proof} By inspecting the definitions, we get an
isomorphism of $B_n$-modules
\set\begin{equation}\label{eq_shift-module}
A_{[\gamma]}\isom\bigoplus_{\delta\in\Sigma}
\{\delta\otimes a\in K[\Delta^\gp_n]~|~\log|a|\geq
f_A(\gamma+\delta)\}
\end{equation}
for every $\Sigma\in\cS_{\!\Delta^\gp_n,\Delta_\R}$ and every
$\gamma\in\Omega(\Delta_\R,\Sigma)\cap\Delta^\gp_\Q$. Hence,
suppose more generally that $\Sigma,\Sigma'$ are two elements of
$\cS_{\!\Delta^\gp_n,\Delta_\R}$ such that $\Sigma\subset\Sigma'$,
and let $\gamma\in\Omega(\Delta_\R,\Sigma)\cap\Delta^\gp_\Q$,
$\gamma'\in\Omega(\Delta_\R,\Sigma')\cap\Delta^\gp_\Q$ be any two
elements; from lemma \ref{lem_liner-subdivide}(ii) we get, for
any $b\in K^+$ with $\log|b|\geq C_A\Vert\gamma'-\gamma\Vert$,
a $B_n$-linear map
$$
\tau_{b,\gamma'-\gamma}:A_{[\gamma]}\to A_{[\gamma']}
$$
which, under the identifications \eqref{eq_shift-module},
corresponds to the map given by the rule :
$\delta\otimes a\mapsto\delta\otimes ba$ for every $\delta\in\Sigma$
and every $a\in K$ such that $\log|a|\geq f(\gamma+\delta)$.
In case $\Sigma=\Sigma'$, also $\tau_{b,\gamma-\gamma'}$ is
well defined, and clearly
$$
\tau_{b,\gamma'-\gamma}\circ\tau_{b,\gamma-\gamma'}=
b^2\cdot\one_{A_{[\gamma']}}
\qquad
\tau_{b,\gamma-\gamma'}\circ\tau_{b,\gamma'-\gamma}=
b^2\cdot\one_{A_{[\gamma]}}
$$
whence the contention.
\end{proof}

\begin{remark}\label{rem_four-gammas}
Let $\gamma_1,\gamma_2,\gamma'_1,\gamma'_2\in\Delta^\gp_\Q$
be four elements, and set $\gamma_3:=\gamma_1+\gamma_2$,
$\gamma'_3:=\gamma'_1+\gamma'_2$; suppose that
$$
\Delta^\gp_n\cap(\Delta_\R-\gamma_i)\subset
\Delta^\gp_n\cap(\Delta_\R-\gamma'_i)
\qquad
\text{for $i=1,2,3$}.
$$
With the notation of the proof of lemma \ref{lem_lot-of-notate},
for every $b_1,b_2\in K^+$ such that
$$
\log|b_i|\geq C_A\Vert\gamma'_i-\gamma_i\Vert
\qquad
(i=1,2)
$$
we obtain a commutative diagram of $B_n$-linear maps
$$
\xymatrix{
A_{[\gamma_1]}\otimes_{B_n}A_{[\gamma_2]} \ar[r]
\ar[d]_{\tau_{b_1,\gamma'_1-\gamma_1}\otimes\tau_{b_2,\gamma'_2-\gamma_2}} &
A_{[\gamma_3]} \ar[d]^{\tau_{b_1b_2,\gamma'_3-\gamma_3}}
\\
A_{[\gamma'_1]}\otimes_{B_n}A_{[\gamma'_2]} \ar[r] & A_{[\gamma'_3]}
}$$
whose horizontal arrows are the restrictions of the multiplication
map of the $B_n$-algebra $B$ : the detailed verification shall be
left as an exercise for the reader.
\end{remark}

\begin{theorem}\label{th_B-sh-is-measurable}
In the situation of \eqref{subsec_extend-to-Q}, the
$K^+$-algebra $B^\sh$ is ind-measurable.
\end{theorem}
\begin{proof} First, notice that, for every $k\geq 0$, the
inclusion map $B_n\to B_{n+k}$ induces a radicial morphism
$B_n\otimes_{K^+}\kappa\to B_{n+k}\otimes_{K^+}\kappa$,
and therefore also an isomorphism of $B_{n+k}$-algebras
\set\begin{equation}\label{eq_radicial-B}
B_n^\sh\otimes_{B_n}B_{n+k}\isom B^\sh_{n+k}
\end{equation}
(\cite[Ch.IV, Prop.18.8.10]{EGA4}). Notice as well, that
$B^\sh_n$ is a measurable $K^+$-algebra (see definition
\ref{def_working-algeb}(ii)), hence we have a well defined
normalized length function $\lambda_n$ on isomorphism classes
of $B_n^\sh$-modules, characterized as in theorem
\ref{th_about_K_0}. Thus, our task is to exhibit a sequence
$(d_n~|~n\in\N)$ of normalizing factors fulfilling conditions
(a) and (b) of definition \ref{def_ind-measur}, for the directed
system $(B_n^\sh~|~n\in\N)$, whose colimit is $B^\sh$.
To this aim, let $M$ be any object of $B^\sh_n\Mod_\cohs$;
in view of \eqref{eq_radicial-B}, we need to estimate
$$
\lambda_{k+n}(B_{k+n}\otimes_{B_n}M)=
\lambda_{k+n}
\left(\bigoplus_{[\gamma]\in\Delta^\gp_{k+n}/\Delta^\gp_n}
A_{[\gamma]}\otimes_{B_n}M\right).
$$
To this aim, let us introduce the function
$$
l_M:\Delta^\gp_\Q\to\R
\qquad
\gamma\mapsto\lambda_n(A_{[\gamma]}\otimes_{B_n}M).
$$

\begin{claim}\label{cl_lipsync}
Let $\Sigma$ be as in lemma \ref{lem_lot-of-notate}, set
$J:=\Ann_{B^\sh_0}(M/\fm M)$, and suppose that $M$ admits
a generating set of cardinality $k$. Then the restriction
of $l_M$ to $\Omega(\Delta_\R,\Sigma)\cap\Delta^\gp_\Q$ is
{\em of Lipschitz type}, {\em i.e.} there exists a real
constant $C'_A>0$ independent of $n$, $\Sigma$ and $M$ such that
$$
|l_M(\gamma)-l_M(\gamma')|\leq
C'_A\cdot k\cdot\length_{B^\sh_n}(B^\sh_n/JB^\sh_n)\cdot
\Vert\gamma-\gamma'\Vert
$$
for every $\gamma,\gamma'\in\Omega(\Delta_\R,\Sigma)\cap\Delta^\gp_\Q$.
\end{claim}
\begin{pfclaim} Let $I\subset B^\sh_0$ be any finitely generated
$\fm_{B_n^\sh}$-primary ideal contained in the annihilator of $M$.
Also, let $N'$ be the cardinality of a finite set of generators
for the $B^\sh_n$-module $M$, and $N$ the constant provided by
lemma \ref{lem_choose-int-N}; clearly
$A_{[\gamma]}\otimes_{B_n}M$ admits a system of generators of
cardinality $\leq NN'$, for every $\gamma\in\Delta_\Q^\gp$.
Furthermore, an argument as in the proof of
\cite[Lemma 2.3.7(iv)]{Ga-Ra} shows that the induced mapping
$$
\cM_N(B_n)\to\cM_{NN'}(B^\sh_n/IB_n^\sh)
\qquad
L\mapsto L\otimes_{B_0}M
$$
is of Lipschitz type; more precisely, it sends the entourage
$E_r$ into $E_{2r}$, for every $r\in\R_{>0}$ (details left to
the reader). Then the claim follows by combining lemmata
\ref{lem_ta-ta} and \ref{lem_lot-of-notate}.
\end{pfclaim}

From proposition \ref{prop_linear-constr-part}(iv), we see
that, for every $\Sigma\in\cS_{\!\Delta^\gp_n,\Delta_\R}$, the
subset $\Omega(\Delta_\R,\Sigma)\cap\Delta^\gp_\Q$ is dense in
$\Omega(\Delta_\R,\Sigma)$; taking into account claim
\ref{cl_lipsync}, it follows that $l_M$ extends (uniquely)
to a function $l_{M,\R}:\Delta^\gp_\R\to\R$, whose restriction
to every $\Omega(\Delta_\R,\Sigma)$ is continuous and still
satisfies the same Lipschitz type estimate. Notice that
$l_{M,\R}$ descends to a function on the torus
$$
\bar l_M:\bT_n:=\Delta^\gp_\R/\Delta^\gp_n\to\R.
$$
Hence, set $r:=\rk_\Z\Delta^\gp_0$, let us fix a basis
$e_1,\dots,e_r$ for the free $\Z$-module $\Delta^\gp_0$,
and define
$$
\Omega_n:=
\left\{
\sum_{i=1}^rt_ie_i~|~-\frac{1}{2p^n}\leq t_i<\frac{1}{2p^n}\ \
\text{for $i=1,\dots,r$}
\right\}
$$
so $\Omega_n\subset\Delta^\gp_\R$ is a {\em fundamental domain\/}
for the lattice $\Delta^\gp_n$; in view of theorem
\ref{th_about_K_0}(ii.b), we reach the following identity :
$$
\lambda_{k+n}(B_{k+n}\otimes_{B_n}M)=
\frac{1}{[\kappa(B^\sh_{k+n}):\kappa(B^\sh_n)]}\cdot
\sum_{\gamma\in\Delta^\gp_{k+n}\cap\Omega_n} l_M(\gamma).
$$
Now, set
\set\begin{equation}\label{eq_norm-factors-B}
d_n:=p^{nr}\cdot[\kappa(B^\sh_n):\kappa(B^\sh_0)]^{-1}
\qquad
\text{for every $n\in\N$}.
\end{equation}
We claim that $(d_n~|~n\in\N)$ is a suitable sequence of
normalizing factors for $B^\sh$. Indeed, recall that
$\Omega(\Delta_\R,\Sigma)$ is linearly constructible for every
$\Sigma\in\cS_{\!\Delta^\gp_0,\Delta_\R}$ (proposition
\ref{prop_linear-constr-part}(iii)), especially, the
bounded function $l_M$ is continuous outside a subset of
$\Delta^\gp_\R$ of measure zero, hence it is integrable
on every bounded measurable subset of $\Delta_\R$.
It follows that $\bar l_M$ is integrable relative to
the invariant measure $d\mu_n$ on $\bT_n$ of total volume
equal to $1$; lastly, a simple inspection yields
$$
\lambda(B\otimes_{B_n}M):=\lim_{k\to+\infty}
d_{k+n}^{-1}\cdot\lambda_{k+n}(B_{k+n}\otimes_{B_n}M)=
d_n^{-1}\int_{\bT_n}\bar l_Md\mu_n
$$
which shows that condition (a) of definition
\ref{def_ind-measur} holds for this choice of normalizing
factors.

In order to check condition (b), set
$\Delta_\R(\rho):=\{v\in\Delta_\R~|~\Vert v\Vert\leq\rho\}$
for every $\rho>0$. Notice that the automorphism of
$\Delta^\gp_\R$ given by multiplication by $p^n$ restricts
to a bijection
$$
\Omega(\Delta_\R,\Delta_n)\isom\Omega(\Delta_\R,\Delta_0)
\qquad
\text{for every $n\in\N$}.
$$
Then, by claim \ref{cl_C_eps-lies}, we see that there exists
$\rho_0$ such that
\set\begin{equation}\label{eq_small-rho}
\Delta_\R(p^{-n}\rho_0)\subset\Omega_n\cap\Omega(\Delta_\R,\Delta_n)
\qquad
\text{for all $n\in\N$}.
\end{equation}
Especially, for every $\rho\leq\rho_0$ we may regard
$\Delta_\R(p^{-n}\rho)$ as a measurable subset of $\bT_n$,
whose measure $\mathrm{Vol}(\rho)$ is strictly positive
and independent of $n$.

\begin{claim}\label{cl_cons-with-I}
Let $I\subset B^\sh_0$ be any finitely generated
$\fm_{B^\sh_0}$-primary ideal. Then there exists
a real constant $C_I>0$ such that
$$
d_n^{-1}\cdot\length_{B^\sh_n}(B_n^\sh/(I+\fm)B^\sh_n)\leq C_I
\qquad
\text{for every $n\in\N$}.
$$
\end{claim}
\begin{pfclaim} On the one hand, we may write $B^\sh_n$
as the direct sum of the $B_0^\sh$-modules
$A_{[\gamma]\cap\Delta_\Q}\otimes_{B_0}B^\sh_0$, for
$\gamma$ ranging over the elements of $\Delta^\gp_n/\Delta^\gp_0$.
There are $p^{rn}$ such direct summands, and each of them
admits a generating system of cardinality $\leq N$, where
$N$ is the constant provided by lemma \ref{lem_choose-int-N}.
It follows that
$$
\length_{B^\sh_0}(B^\sh_n/(I+\fm)B_n^\sh)\leq
Np^{rn}\cdot\length_{B_0^\sh}(B^\sh_0/(I+\fm B_0^\sh)).
$$
On the other hand, say that
$l:=\length_{B^\sh_n}(B_n^\sh/(I+\fm)B^\sh_n)$; this means
that may find a filtration of $B_n^\sh/(I+\fm)B^\sh_n$ of
length $l$, consisting of $B_n^\sh$-submodules, whose
graded subquotients are all isomorphic to $\kappa(B^\sh_n)$.
Given such a filtration, we easily obtain a filtration
of $B_n^\sh/(I+\fm)B^\sh_n$ of length
$l\cdot[\kappa(B^\sh_n):\kappa(B^\sh_0)]$, consisting of
$B_0^\sh$-submodules, whose graded subquotients are all
isomorphic to $\kappa(B^\sh_0)$. Therefore
$$
l\cdot[\kappa(B^\sh_n):\kappa(B^\sh_0)]=
\length_{B^\sh_0}(B^\sh_n/(I+\fm)B_n^\sh)
$$
whence the claim.
\end{pfclaim}

Now, fix $\eps>0$, and suppose there is given a finitely
generated $\fm_{B_0^\sh}$-primary ideal $I\subset B^\sh_0$,
and a surjection of finitely presented $B^\sh_n/IB^\sh_n$-modules
$M\to M'$, generated by $k$ elements, such that
$$
d_n^{-1}\cdot(\lambda_n(M)-\lambda_n(M'))>\eps.
$$
Set $C(I,k):=2kC'_AC_I$, where $C'_A$ and $C_I$ are as in
claims \ref{cl_lipsync} and \ref{cl_cons-with-I}; since
$$
l_M(0)-l_{M'}(0)=\lambda_n(M)-\lambda_n(M')
$$
we may estimate :
$$
\begin{aligned}
|\lambda(B^\sh\otimes_{B^\sh_n}M)-
\lambda(B^\sh\otimes_{B^\sh_n}M')|
\geq\: & d_n^{-1}\cdot
\int_{\Delta_\R(p^{-n}\rho)}(\bar l_M-\bar l_{M'})d\mu_n \\
\geq\: &
\mathrm{Vol}(\rho)\cdot(d_n^{-1}(l_M(0)-l_{M'}(0))-
                          C(I,k)\cdot p^{-n}\rho) \\
\geq\: &
\mathrm{Vol}(\rho)\cdot(\eps-C(I,k)\cdot\rho)
\end{aligned}
$$
for every $\rho\leq\rho_0$. Therefore, if we set
$$
\rho(k,\eps,I):=\min\{\rho_0,2^{-1}C(I,k)^{-1}\eps\}
$$
it is easily seen that the sought condition (b) holds with
$\delta(k,\eps,I):=2^{-1}\cdot\mathrm{Vol}(\rho(k,\eps,I))\cdot\eps$.
\end{proof}

\subsection{Almost purity : the log regular case}
\label{sec_regular-pure}
In this section, we prove an almost purity theorem for certain
towers of regular log schemes.

\sset\subsubsection{}\label{subsec_C_n}
Let $\underline M{}_0$ be a log structure on the Zariski site
of a local scheme $X_0$, such that $(X_0,\underline M{}_0)$ is
regular at the closed point $x_0\in X_0$, and say that
$X_0=\Spec\,B_0$ for a local ring $B_0$ which is necessarily
noetherian, normal and Cohen-Macaulay (corollary
\ref{cor_normal-and-CM}). Let also $\beta_0:P\to B_0$ be
a chart for $\underline M{}_0$ that is sharp at $x_0$.
Especially, $P$ is a fine and saturated monoid, and
$A:=B_0/\fm_P B_0$ is a regular local ring; we denote by $\fm_A$
(resp. $\fm_{B_0}$) the maximal ideal of $A$ (resp. of $B_0$).
We assume furthermore that :
\begin{enumerate}
\alphaenu
\item
The characteristic $p$ of the residue field $\kappa(x_0)$ of $A$
is positive.
\item
The Frobenius endomorphism $\Phi_{B_0}$ of $B_0/pB_0$ is a finite
ring homomorphism.
\end{enumerate} 
Notice that (b) implies that $B_0/pB_0$ is excellent (theorem
\ref{th_Kunz-exc}(i)) and also that $\Omega^1_{A/\Z}\otimes_A\kappa(x_0)$
is a finite dimensional $\kappa(x_0)$-vector space.
Moreover, $(X_0,\underline M{}_0)$ is a regular log
scheme, by theorem \ref{th_reg-generizes}.

\begin{example}\label{ex_basic-pairs}
Let $(B,\beta)$ be a pair consisting of
\begin{enumerate}
\alphaenu
\item
a local ring $B$ whose residue field $\kappa_B$ has
positive characteristic
\item
a local morphism of monoids $\beta:P\to B$ from a fine,
sharp and saturated monoid $P$.
\end{enumerate}
Denote by $\fm_B$ the maximal ideal of $B$, by $\underline M$
the logarithmic structure on the Zariski site of $X:=\Spec\,B$
induced by $\beta$, and suppose that $(X,\underline M)$ is
regular at the closed point of $X$.

(i)\ \
If moreover, $B$ is complete for the $\fm_B$-adic topology
and $\kappa_B$ is perfect, then the Frobenius endomorphism
of $B/pB$ is a finite map. Indeed, in this case, for a
suitable $d\in\N$ we have a surjective ring homomorphism
$$
R[[P\times\N^{\oplus d}]]\to B
$$
where $R$ is either $\kappa_B$, or a discrete valuation
ring with an isomorphism $R/pR\isom\kappa_B$ (remark
\ref{rem_more-precisely}). The assertion is an immediate
consequence.

(ii)\ \
It follows that for any such pair $(B,\beta)$ we may
find a local ring $B_0$ and a local and flat ring
homomorphism $f:B\to B_0$ such that the log scheme
$$
(X_0,\underline M{}_0):=\Spec\,B_0\times_X(X,\underline M)
$$
and its chart $\beta_0:=\beta\circ f:P\to B_0$ fulfill the
conditions of \eqref{subsec_C_n}. Indeed, according to
\cite[Ch.0, Prop.10.3.1]{EGAIII} there exists a
flat and local $B$-algebra $C$ such that
$\kappa_B\otimes_BC$ is a perfect field, and we take
$B_0$ to be the completion of $C$. By proposition
\ref{prop_we-got-a-situation}(ii), the resulting log scheme
$(X_0,\underline M{}_0)$ is regular at the closed point
of $X_0$, so the assertion follows from (i).
\end{example}

\sset\subsubsection{}\label{subsec_more-notation}
In the situation of \eqref{subsec_C_n}, set
$$
P^{(n)}:=\{\gamma\in P_\Q~|~\gamma^{p^n}\in P\}
\qquad
\text{for every $n\in\N$}.
$$
The $p^n$-Frobenius map of $P^{(n)}$ identifies $P^{(n)}$ with its
submonoid $P$; in other words, the inclusion map $P\to P^{(n)}$
is naturally identified with the $p^n$-Frobenius endomorphism of
$P$. Fix a sequence $(f_1,\dots,f_r)$ of elements of $B_0$ whose
image in $A$ is maximal in the sense of remark
\ref{rem_afterthought}(iii), in which case we say that
$(f_1,\dots,f_r)$ is {\em maximal in $B_0$}. Notice that
$\dim_{\kappa(x_0)}\fm_A/\fm_A^2=\dim A$, since $A$ is regular; it
follows that
$$
r=\dim A+\dim_{\kappa(x_0)}\Omega^1_{\kappa(x_0)/\Z}
$$
by virtue of the short exact sequence \eqref{eq_finally-exact},
if $p\in\fm_A^2$, and otherwise, by virtue of proposition
\ref{prop_lift-lemma}.
Denote also by $\bar f_i\in A$ the image of $f_i$, for every
$i=1,\dots,r$. According to corollary \ref{cor_regu-criterion},
the ring
$$
A_n:=A[T_1,\dots,T_r]/(T_1^{p^n}-\bar f_1,\dots,T_r^{p^n}-\bar f_r)
$$
is regular. For every $n\in\N$, we set :
$$
B'_n:=P^{(n)}\otimes_PB_0
\qquad
B''_n:=B_0[T_1,\dots,T_r]/(T_1^{p^n}-f_1,\dots,T_r^{p^n}-f_r)
\qquad
B_n:=B'_n\otimes_{B_0}B''_n.
$$

\begin{lemma}\label{lem_take-this}
The induced maps
$$
B_n\to B_{n+1}
\qquad\text{and}\qquad
B_n/pB_n\to B_{n+1}/pB_{n+1}
$$
are injective, for every $n\in\N$. 
\end{lemma}
\begin{proof} The natural map $B_n\to B_{n+1}$ factors as the
composition
$$
B'_n\otimes_{B_0}B''_n\to
B'_{n+1}\otimes_{B_0}B''_n\to B'_{n+1}\otimes_{B_0}B''_{n+1}
\qquad
\text{for every $n\in\N$}.
$$
Notice that $B_0=B''_0$, and $B''_{n+1}$ is a free $B''_n$-module (of
rank $p^r$) for every $n\in\N$; we are then reduced to checking
that the maps $B'_n\to B'_{n+1}$ and
$B'_n/pB'_n\to B'_{n+1}/pB'_{n+1}$ are injective for every $n\in\N$.
However, set $G:=P^{(n+1)}/P^{(n)}$, and notice that $P^{(n+1)}$
is a $G$-graded monoid, with $P^{(n+1)}_0=P^{(n)}$, hence
$B'_{n+1}$ is a $G$-graded $B'_n$-algebra with $(B'_{n+1})_0=B'_n$
for every $n\in\N$. The assertion follows.
\end{proof}

\sset\subsubsection{}\label{subsec_another-Frobenius}
In view of lemma \ref{lem_take-this}, we may set
$$
B'':=\bigcup_{n\in\N}B''_n
\qquad
B:=\bigcup_{n\in\N}B_n
\qquad
P^{(\infty)}:=\bigcup_{n\in\N}P^{(n)}
$$
and clearly
$$
B=P^{(\infty)}\otimes_PB''
\qquad
\text{for every $n\in\N$}
$$
from which we see that $B$ is naturally a
$P^{(\infty)}/P^{(n)}$-graded $B_n$-algebra, for every
$n\in\N$. Also, the induced morphism
$\Spec\,B_{n+1}/pB_{n+1}\to\Spec\,B_n/pB_n$ is radicial
and surjective, so $\Spec\,B_n/pB_n$ is a local scheme
for every $n\in\N$; on the other hand, the map $B_0\to B_n$
is finite, so every point of $\Spec\,B_n$ specializes
to a point of $\Spec\,B_n/pB_n$. We conclude that $B_n$
is a local ring, and we denote by $\fm_{B_n}$ its maximal
ideal, for every $n\in\N$.
Let $\underline M{}_n$ be the log structure on the Zariski
site of $X_n:=\Spec\,B_n$ deduced from the natural map
$\beta_n:P^{(n)}\to B_n$; notice that
$B_n/\fm_{P^{(n)}}B_n=A_n$ is a regular local ring of
dimension equal to $\dim A$. Since we have as well
$\dim P^{(n)}=\dim P$, it follows that $(X_n,\underline M{}_n)$
is regular at the closed point $x_n\in X_n$. Then theorem
\ref{th_reg-generizes} shows that $(X_n,\underline M{}_n)$
is a regular log scheme. Thus, we have obtained a tower
of finite morphisms of regular log schemes
\set\begin{equation}\label{eq_max-tower}
\cdots\to(X_{n+1},\underline M{}_{n+1})\to(X_n,\underline M{}_n)
\to\cdots\to(X_0,\underline M{}_0)
\end{equation}
which we call the {\em maximal tower\/} associated with the
chart $\beta_0:P\to B_0$ and the maximal sequence $(f_1,\dots,f_r)$.
The limit of the tower \eqref{eq_max-tower} is a log scheme
$(X,\underline M)$ whose log structure admits a chart
$P^{(\infty)}\to B$ which is sharp at the closed point.

\begin{remark}\label{rem_shift-sequence}
Keep the notation of \eqref{subsec_another-Frobenius}, and let
$s\in\N$ be any integer; from remark \ref{rem_afterthought}(iii)
it is easily seen that the sequence
$((X_{n+s},\underline M{}_{n+s})~|~n\in\N)$, obtained by removing
from \eqref{eq_max-tower} the first $s$ terms, is the maximal
tower associated with the chart $\beta_n$ and the maximal sequence
$(f_1^{1/p^s},\dots,f_r^{1/p^s})$.
\end{remark}

\sset\subsubsection{}\label{subsec_local-max}
Let now $y\in(X_0,\underline M{}_0)_\tr$ be any point, and
$n\in\N$ any integer. The chart $\beta_0$ extends uniquely
to a morphism of monoids $\beta^\gp_0:P^\gp\to B_{0,y}:=\cO_{\!X_0,y}$,
and we have a natural isomorphism
$$
B_n\otimes_{B_0}B_{0,y}\isom
P^{(n)\gp}\otimes_{P^\gp}B_{0,y}\otimes_{B_0}B''_n.
$$
Especially, $B_n\otimes_{B_0}B_{0,y}$ is a free $B_{0,y}$-module
of rank $p^{nd}$, where
$$
d:=\dim P+r=
\dim P+\dim A+\dim_{\kappa(x_0)}\Omega^1_{\kappa(x_0)/\Z}.
$$
However, since $(X_0,\underline M{}_0)$ is regular, we have
$\dim A+\dim P=\dim B_0$. Summing up, we find
$$
d=\dim B_0+\dim_{\kappa(x_0)}\Omega^1_{\kappa(x_0)/\Z}.
$$
Next, let $z\in\Spec\,B_0/pB_0\subset X_0$ be any point, and
$\fp_z\subset B_0$ the corresponding prime ideal; as in the
foregoing, the chart $\beta_0$ extends uniquely to a morphism
$\beta_\fp:P_\fp\to B_{0,z}:=\cO_{\!X_0,z}$, where $\fp:=\beta_0^{-1}\fp_z$.
As already remarked, the point $z$ lifts uniquely to a point
$z_n\in X_n$, and on the other hand, there exists a unique
prime ideal $\fp^{(n)}\subset P^{(n)}$ containing $\fp$,
and the inclusion map $j_\fp:P_\fp\to P^{(n)}_{\fp^{(n)}}$ is
naturally identified with the $p^n$-Frobenius endomorphism
of $P_\fp$.
By lemma \ref{lem_decomp-sats}, there exists an isomorphism
of monoids $P_\fp\isom G\times Q$, with $G:=P_\fp^\times$
and $Q:=P_\fp^\sharp$; we may then find a corresponding
decomposition $P^{(n)}_{\fp^{(n)}}=G^{(n)}\times Q^{(n)}$
that identifies $j_\fp$ with the product of maps of monoids
$G\to G^{(n)}$ and $Q\to Q^{(n)}$. Summing up, there follows
an isomorphism of $B_0$-algebras
$$
\cO_{\!X_n,z_n}\isom
P^{(n)}_{\fp^{(n)}}\otimes_{P_\fp}B_{0,z}\otimes_{B_0}B''_n\isom
(Q^{(n)}\otimes_QB_{0,z})\otimes_{B_{0,z}}(G^{(n)}\otimes_GB_{0,z})
\otimes_{B_0}B''_n.
$$
Fix a basis $g_1,\dots,g_s$ of the free $\Z$-module $G$, and
set $f_{r+i}:=\beta_\fp(g_i)$ for $i=1,\cdots,s$; clearly
$$
(G^{(n)}\otimes_GB_{0,z})\otimes_{B_0}B''_n=
B_{0,z}[T_1,\dots,T_{r+s}]/(T_1^{p^n}-f_1,\dots,T_{r+s}^{p^n}-f_{r+s}).
$$
On the other hand, set $A_z:=B_{0,z}/\fm_QB_{0,z}$; we have
$$
\begin{aligned}
\dim P=\: & \dim Q+s & &
\qquad\text{(corollary \ref{cor_consequent}(i))} \\
\dim B_{0,z}=\: &\dim Q+\dim A_z & &
\qquad\text{(since $(X_0,\underline M{}_0)$ is regular)} \\
\dim B_{0,z}=\: & d-\dim_{\kappa(z)}\Omega^1_{\kappa(z)/\Z}
& & \qquad\text{(proposition \ref{prop_jump-of-the-way})}.
\end{aligned}
$$
Therefore $r+s=\dim A_z+\dim_{\kappa(z)}\Omega^1_{\kappa(z)/\Z}$.
However, $A_{z_n}:=\cO_{\!X_n,z_n}/\fm_{Q^{(n)}}\cO_{\!X_n,z_n}$
is a regular local ring, since $(X_n,\underline M{}_n)$
is regular, and by inspecting the construction we see that
$$
A_{z_n}=
A_z[T_1,\dots,T_{r+s}]/(T_1^{p^n}-f_1,\dots,T_{r+s}^{p^n}-f_{r+s})
$$
hence the image of the system $(f_1,\dots,f_{r+s})$ yields a
basis of either $\Omega^1_{A_z/\Z}\otimes_{A_z}\kappa(z)$ or
$\bOmega_{A_z}$, depending on whether or not $p\in\fm^2_{A_z}$
(corollary \ref{cor_regu-criterion}). In conclusion, we see
that the induced sequence of morphisms of log schemes
$$
\cdots\to(X_{n+1}(z_n),\underline M{}_{n+1}(z_n))\to
(X_n(z_n),\underline M{}_n(z_n))
\to\cdots\to(X_0(z),\underline M_0(z))
$$
(notation of \eqref{subsec_acyclic-log}) is the maximal tower
associated with the induced chart $Q\to\cO_{\!X_0,z}$ and the maximal
sequence $(f_1,\dots,f_{r+s})$.

\begin{remark}\label{rem_reduce-to-sh}
(i)\ \
Let $\bar x_0$ be any geometric point of $X_0$, localized
at $x_0$, and denote by $B_0^\sh$ the strict henselization
of $B_0$ at the point $\bar x_0$. In light of
\cite[Th.3.5.13(ii)]{Ga-Ra}, we see that the log scheme
$(X_0(\bar x_0),\underline M{}_0(\bar x_0))$ still fulfills
the conditions of \eqref{subsec_C_n}. Moreover, the image
in $B^\sh_0$ of the maximal sequence $(f_1,\dots,f_n)$ is
still maximal (remark \ref{rem_afterthought}(iv)). It
follows easily that the induced tower of regular log
schemes
$$
\cdots\to X_0(\bar x_0)\times_{X_0}(X_{n+1},\underline M{}_{n+1})
\to X_0(\bar x_0)\times_{X_0}(X_n,\underline M{}_n)\to\cdots\to
(X_0(\bar x_0),\underline M_0(\bar x_0))
$$
is (naturally isomorphic to) the maximal tower
associated with the induced chart $P_0\to B^\sh_0$
and the maximal sequence $(f_1,\dots,f_n)$ of $B^\sh_0$.

(ii)\ \
In the same vein, let $I\subset B_0$ be any proper ideal,
$B^\wedge_{0,I}$ the $I$-adic completion of the local
ring $B_0$, and set $X_0^\wedge:=\Spec\,B^\wedge_{0,I}$;
according to claim \ref{cl_stability}(i), the log
scheme $(X_0^\wedge,\underline M{}_0^\wedge):=
X_0^\wedge\times_{X_0}(X_0,\underline M{}_0)$ still
fulfills the conditions of \eqref{subsec_C_n}, and
the image in $B^\wedge_{0,I}$ of the maximal sequence
$(f_1,\dots,f_n)$ of $B_0$ is still maximal
(remark \ref{rem_afterthought}(iv)). It follows
that the induced tower of log schemes
$$
\cdots\to X^\wedge_0\times_{X_0}(X_{n+1},\underline M{}_{n+1})
\to X_0^\wedge\times_{X_0}(X_n,\underline M{}_n)\to\cdots\to
(X_0^\wedge,\underline M{}_0^\wedge)
$$
is  (naturally isomorphic to) the maximal tower
associated with the induced chart $P_0\to B^\wedge_{0,I}$
and the maximal sequence $(f_1,\dots,f_n)$ of
$B^\wedge_{0,I}$.
\end{remark}

\sset\subsubsection{}
Let $\Phi_{B_n}:B_n/pB_n\to B_n/pB_n$ be the Frobenius
endomorphism of $B_n/pB_n$; taking into account lemma
\ref{lem_take-this}, we see that $\Phi_{B_{n+1}}$ factors
through a unique ring homomorphism
$$
\Psi_{B_{n+1}}:B_{n+1}/pB_{n+1}\to B_n/pB_n
\qquad
\text{for every $n\in\N$}.
$$

\begin{lemma}\label{lem_Frob-onto}
The map $\Psi_{B_{n+1}}$ is surjective for every $n\in\N$.
\end{lemma}
\begin{proof} Let us start out with the following general :

\begin{claim}\label{cl_radicial-etale}
Let $\phi:R\to S$ be an injective, finite and radicial ring
homomorphism. Then $\phi$ is an isomorphism if and only if
$\Omega^1_{S/R}=0$.
\end{claim}
\begin{pfclaim} We may assume that $\Omega^1_{S/R}=0$, and
we show that $\phi$ is an isomorphism. To this aim, it
suffices to show that
$$
\phi\otimes_RR/\fm:R':=R/\fm\to S':=S/\fm S
$$
is an isomorphism for every maximal ideal $\fm\subset R$.
However, $\Omega^1_{S'/R'}=\Omega^1_{S/R}\otimes_RR'=0$,
so we may replace $R$ by $R'$ and $S$ by $S'$, and assume
from start that $R$ is a field. In this case, $S$ is a
local unramified $R$-algebra, so $S$ must be a finite
separable field extension of $R$, by
\cite[Ch.IV, Cor.17.4.2]{EGA4}. But $S$ is also a radicial
extension of $R$, so $S=R$.
\end{pfclaim}

\begin{claim}\label{cl_Anatole}
Let $p>0$ be a prime integer, $R$ a local $\F_p$-algebra
whose Frobenius endomorphism $\Phi_R$ is a finite map,
$k_R$ the residue field of $R$, and $\Sigma\subset R$ any
subset. Set $R_0:=R/\Sigma R$, and let $(g_1,\dots,g_n)$
be a finite sequence of elements of $R$ such that
$dg_1,\dots,dg_n$ is a system of generators for the
$k_R$-vector space $\Omega^1_{R_0/\Z}\otimes_Rk_R$. Then
$R=R^p[\Sigma,g_1,\dots,g_n]$.
\end{claim}
\begin{pfclaim} Set $S:=R^p[\Sigma,g_1,\dots,g_n]$. The
inclusion map $S\to R$ is clearly radicial, and it is finite,
since $\Phi_R$ is finite; hence claim \ref{cl_radicial-etale}
reduces to checking that $\Omega^1_{R/S}$ vanishes.
By Nakayama's lemma, it then suffices to show that
$\Omega^1_{R/S}\otimes_Rk_R=0$. However, let $M$ (resp. $N$)
be the $R$-submodule of $\Omega:=\Omega^1_{R/\Z}$ generated
by $\{da~|~a\in\Sigma\}$ (resp. by $\{dg_1,\dots,dg_n\}$); then
$\Omega^1_{R/S}=\Omega/(M+N)$ (\cite[Ch.0, Th.20.5.7(i)]{EGAIV})
and on the other hand, the induced map
$(\Omega/M)\otimes_RR_0\to\Omega^1_{R_0/\Z}$ is an isomorphism
(\cite[Ch.0, Th.20.5.12(i)]{EGAIV}). We deduce a right exact
sequence of $R$-modules
$$
N\xrightarrow{\ j\ }\Omega^1_{R_0/\Z}\otimes_Rk_R\to
\Omega^1_{R/S}\otimes_Rk_R\to 0.
$$
But our choice of the sequence $(g_1,\dots,g_n)$ implies that
$j$ is surjective, whence the contention.
\end{pfclaim}

Finally, notice that the image of $\Psi_{B_{n+1}}$ is the
subring $(B_n/pB_n)^p[\beta_n(P^{(n)}),f^{1/p^n}_1,\dots,f^{1/p^n}_r]$.
For every $n>0$, consider the exact sequence
$$
\Omega^1_{A/\Z}\otimes_A\kappa(x_n)\xrightarrow{\ j\ }
\Omega^1_{A_n/\Z}\otimes_{A_n}\kappa(x_n)\to
\Omega^1_{A_n/A}\otimes_{A_n}\kappa(x_n)\to 0
$$
(\cite[Ch.0, Th.20.5.7(i)]{EGAIV}); the image of $j$ is generated
by the image of the generating system $df_\bullet:=(df_1,\dots,df_r)$
of the $\kappa(x_0)$-vector space $\Omega^1_{A/\Z}\otimes_A\kappa(x_0)$.
However, since $f_i$ admits a $p$-th root in $A_n$ for every
$i=1,\dots,r$, it is easily seen that the image of $df_\bullet$
vanishes in $\Omega^1_{A_n/\Z}$. Hence
$\Omega^1_{A_n/\Z}\otimes_{A_n}\kappa(x_n)=
\Omega^1_{A_n/A}\otimes_{A_n}\kappa(x_n)$ is generated by the
image of the system $(df^{1/p^n}_1,\dots,df^{1/p^n}_r)$.
Therefore, to prove the lemma, it suffices to apply claim
\ref{cl_Anatole} to $R:=B_n/pB_n$, $\Sigma:=\beta_n(\fm_{P^{(n)}})$
and the sequence $(f^{1/p^n}_1,\dots,f^{1/p^n}_r)$.
\end{proof}

\begin{theorem}\label{th_not-an-F_p}
With the notation of \eqref{subsec_another-Frobenius}, we have :
\begin{enumerate}
\item
If $B_0$ is an $\F_p$-algebra, $\Psi_{B_n}$ is an isomorphism
for every $n\in\N$.
\item
If $B_0$ is not an $\F_p$-algebra, then there exist $\pi\in B_1$
and $u\in B^\times_1$ such that $\pi^p=pu$.
\item
For every $\pi\in B_1$ as in {\em (ii)} and every integer $n>0$,
we have $\Ker\,\Psi_{B_n}=\pi B_n/pB_n$.
\item
Let $\cT_{p,B}$ be the $p$-adic topology of $B$. Then
$(B,\cT_{p,B})$ is a formal perfectoid ring.
\end{enumerate}
\end{theorem}
\begin{proof} (i) is immediate from lemma \ref{lem_Frob-onto},
since in that case $\Phi_{B_n}$ is the Frobenius endomorphism
of $B_n$, so it is injective, since the latter is a domain.

(iii): Suppose we have found $\pi$ as in (ii), and let
$x\in B_n$ whose image in $B_n/pB_n$ lies in
$\Ker\,\Psi_{B_n}$; this means that $x^p=py$ holds in
$B_n$ for some $y\in B_n$; hence $(x/\pi)^p\in B_n$, and since
$B_n$ is a normal ring (corollary \ref{cor_normal-and-CM}),
we deduce that $x/\pi\in B_n$, {\em i.e.} $x\in\pi B_n$.

(ii): Assume first that $p\in\fm_A^2$. Then we may write
$$
p=\sum^t_{i=1}b_ib'_i+\sum_{j=1}^sb''_j\cdot\beta_0(x_j)
$$
for certain $b_1,b'_1,\dots,b_t,b'_t\in\fm_{B_0}$,
$b''_1,\dots,b''_s\in B_0$, and $x_1,\dots,x_s\in\fm_P$.
In light of lemma \ref{lem_Frob-onto}, we may then write
$$
b_i=c_i^p+pd_i
\quad
b'_i=c'{}^p_i+pd'_i
\qquad
\text{where $c_i,c'_i,d_i,d'_i\in\fm_{B_1}$ for every $i=1,\dots,t$}
$$
and likewise, $b_j''=c''{}^p_j+pd''_j$ where $c''_j\in B_1$
for $j=1,\dots,s$.
Moreover, by construction, $\beta_0$ extends to a morphism
of monoids $\beta_1:P^{(1)}\to B_1$, and we may write
$x_j=y_j^p$ with $y_j\in\fm_{P^{(1)}}$ for $j=1,\dots,s$.
A simple computation then yields
\set\begin{equation}\label{eq_need-a-break}
p\cdot(1+e)=\sum_{i=1}^tc_i^pc'{}^p_i+
\sum^s_{j=1}c''{}^p_j\cdot\beta_1(y_j)^p
\qquad
\text{for some $e\in\fm_{B_1}$}.
\end{equation}
However, the right-hand side of \eqref{eq_need-a-break} can
be written in the form $g^p+ph$ for some $g,h\in\fm_{B_1}$
(details left to the reader); clearly $1+e-h\in B_1^\times$,
whence the contention, in this case.

Next, suppose that $p\notin\fm_A^2$. In this case, recall
that
$$
\dim_{\kappa(x_0)^{1/p}}\bOmega_A=
1+\dim_{\kappa(x_0)}\Omega^1_{A/\Z}\otimes_A\kappa(x_0).
$$
Therefore, after reordering the sequence $(f_1,\dots,f_r)$,
we may assume that $df_1,\dots,df_{r-1}$ is a basis of
the $\kappa(x_0)$-vector space $\Omega^1_{A/\Z}\otimes_A\kappa(x_0)$.
Set
$$
B':=(P^{(1)}\otimes_PB_0)
[T_1,\dots,T_{r-1}]/(T_1^p-f_1,\dots,T_{r-1}^p-f_{r-1}).
$$
Clearly $B_1$ is a faithfully flat $B'$-algebra, hence
$B'/pB'$ is a $B_0/pB_0$-subalgebra of $B_1/pB_1$, so the natural
map $B_0/pB_0\to B'/pB'$ is injective (lemma \ref{lem_take-this}),
and just as in \eqref{subsec_another-Frobenius}, we deduce
that the Frobenius endomorphism of $B'/pB'$ factors through
a ring homomorphism $\Psi_{B'}:B'/pB'\to B_0/pB_0$, and arguing
as in the foregoing, we also see that $B'$ is a local ring.
Moreover, by applying claim \ref{cl_Anatole} with $R:=B'/pB'$,
$\Sigma:=\beta_0(\fm_P)$ and the sequence $(f_1,\dots,f_{r-1})$
we conclude -- as in the proof of lemma \ref{lem_Frob-onto} --
that $\Psi_{B'}$ is surjective. Hence, denote by $\fm_{B'}$
the maximal ideal of $B'$; it follows that there exist
$g\in\fm_{B'}$ and $h\in B'$ such that $f_r=g^p+ph$ in $B'$.
Set
$$
A':=B'/\fm_{P^{(1)}}B'=
A[T_1,\dots,T_{r-1}]/(T^p_1-f_1,\dots,T^p_{r-1}-f_{r-1}).
$$
Then $A'$ is a regular local ring, by corollary
\ref{cor_regu-criterion}, applied to the sequence
$(\bar f_1,\dots,\bar f_{r-1})$ consisting of the images of
the elements $f_i$ in $A$. By the same criterion -- applied
to the similar sequence $(\bar f_1,\dots,\bar f_r)$ -- we
see that $A'[T]/(T^p-f_r)$ is regular as well. So once
again the same corollary -- applied to the element $\bar f_r$
of $A'$ -- says that the element $\bd(f_r)$ of $\bOmega_{A'}$
does not vanish. However
$$
\bd(f_r)=\bd(g^p)+\bd(ph)=pg^{p-1}\bd(g)+h\bd(p)+p\bd(h)=h\bd(p)
\qquad
\text{in $\bOmega_{A'}$}
$$
(see \eqref{subsec_bar-bd}). We conclude that $h\in B'{}^\times$.
Lastly, notice that $f_r$ admits a $p$-th root $f_r^{1/p}$ in
$B_1$, hence $ph=f_r-g^p=(f_r^{1/p}-g)^p+pe$ in $B_1$ for some
$e\in\fm_{B_1}$ (details left to the reader). Since
$h-e\in B_1^\times$, we are done.

(iv): If $B_0$ is an $\F_p$-algebra, then $\cT_{p,B}$ is
the discrete topology, and $B$ is a perfect $\F_p$-algebra,
whence the assertion in this case. Suppose then that $B_0$
is not an $\F_p$-algebra, in which case the topology of
the completion $B^\wedge$ of $(B,\cT_{p,B})$ is still
$p$-adic (remark \ref{rem_completion-of-topring}(iv)).
Pick $\pi\in B_1$ as in (ii), set $I:=\pi B^\wedge$, and
notice that $\pi$ is a regular element of $B^\wedge$, by
proposition \ref{prop_replaces-Mat-Th.8.1}(i) : the details
shall be left to the reader.
In light of (iii), it is easily seen that $B^\wedge$ is
a P-ring, $I$ is a special ideal of definition, and the
Frobenius endomorphism of $B^\wedge:pB^\wedge=B/pB$ induces
an isomorphism $B^\wedge/I\isom B^\wedge/I^{(p)}$. The assertion
then follows from theorem \ref{th_regular-seq-criterion}.
\end{proof}

\sset\subsubsection{}\label{subsec_beta-is-chart}
Let $\beta:P^{(\infty)}\to B$ be the chart of the log structure
$\underline M$ on $X$ from \eqref{subsec_another-Frobenius},
and notice that the inclusion map $P^{(n)}\to P^{(\infty)}$
induces bijections $\Spec\,P^{(\infty)}\isom\Spec\,P^{(n)}$ for
every $n\in\N$ (lemma \ref{lem_Kummer-fans}(i)); especially,
$\Spec\,P^{(\infty)}$ is a finite set. Let $I\subset B$
be any non-zero ideal; we say that $I$ is a {\em branch ideal},
if there exist radical ideals $J\subset B$ and
$\fr\subset P^{(\infty)}$ such that
\set\begin{equation}\label{eq_branch-ideal}
p\in J
\qquad\text{and}\qquad
I=J\cap\fr B.
\end{equation}

\begin{remark}\label{rem_branch-ideal}
(i)\ \
Notice that, in case $B_0$ is an $\F_p$-algebra, a
branch ideal is just any non-zero radical ideal of $B$.

(ii)\ \
Let $I\subset B$ be a branch ideal, and pick radical
ideals $J\subset B$ and $\fr\subset P^{(\infty)}$ such that
\eqref{eq_branch-ideal} holds.
Set $\fr^{(n)}:=\fr\cap P^{(n)}$ for every $n\in\N$; from corollary
\ref{cor_same-height}(ii) we know that $\fr^{(n)}B_n$
is a radical ideal of $B_n$, hence $\fr B$ is a radical
ideal of $B$, and then the same holds for $I$.

(iii)\ \
In the situation of (ii), let $\fZ\subset\Spec\,P^{(\infty)}$
be a subset such that $\fr=\bigcap_{\fp\in\fZ}\fp$ (lemma
\ref{lem_radical}), and set $\fp^{(n)}:=\fp\cap P^{(n)}$ for
every $\fp\in\fZ$ and every $n\in\N$, so
$\fr^{(n)}=\bigcap_{\fp\in\fZ}\fp^{(n)}$ for every $n\in\N$.
In view of proposition \ref{prop_second-crit} and lemmata
\ref{lem_invariance-of-flatness} and \ref{lem_intersect-ideals},
we see that $\fr^{(n)}B=\bigcap_{\fp\in\fZ}\fp^{(n)}B_n$,
and therefore
\set\begin{equation}\label{eq_radical-identity}
\fr B=\bigcap_{\fp\in\fZ}\fp B.
\end{equation}

(iv)\ \
Furthermore, any branch ideal $I\subset B$ is the radical of
$I_0B$, for some ideal $I_0\subset B_0$.
Indeed, pick $J$ and $\fr$ such that \eqref{eq_branch-ideal}
holds; since the projection $\Spec\,B/pB\to(X_0)_{/p}$ is
radicial and surjective, it is clear that $J$ is the radical
of $J_0B$ for some ideal $J_0\subset B_0$, and on the other
hand, $\fr B$ is the radical of $\fr^{(0)}B$.
\end{remark}

\begin{proposition}\label{prop_branch-setup}
With the notation of \eqref{subsec_beta-is-chart}, let
$I\subset B$ be any branch ideal. We have :
\begin{enumerate}
\item
$I^2=I$, and $I$ fulfills condition $(\bB)$ of\/ \cite[\S2.1.6]{Ga-Ra}.
\item
Let $I_1,I_2\subset B$ be any two branch ideals. Then
$I_1\cap I_2$ and $I_1+I_2$ are also branch ideals, and
$I_1I_2=I_1\cap I_2$.
\end{enumerate}
\end{proposition}
\begin{proof}(i): By assumption, $I=J\cap\fr B$ for radical
ideals $J\subset B$ and $\fr\subset P^{(\infty)}$ with $p\in J$.

\begin{claim}\label{cl_two-squares}
$(\fr B)^2=\fr B$ and $J^2=J$, and the ideals $\fr B$ and
$J$ fulfill condition $(\bB)$.
\end{claim}
\begin{pfclaim} Since the $p$-Frobenius endomorphism
of $P^{(\infty)}$ is an automorphism, the first stated
identity is clear. Next, set $\bar J:=J/pB$.
By lemma \ref{lem_Frob-onto}, the Frobenius endomorphism of
$B/pB$ is surjective; since $\bar J$ is a radical ideal, we
deduce that $\bar J{}^2=\bar J$. On the other hand, from
theorem \ref{th_not-an-F_p}(ii) we see that $p\in J^2$,
from which also the second identity follows easily.

Next, it is clear that $\fr B$ fulfills condition $(\bB)$,
and for $J$ we apply \cite[Claim 2.1.9]{Ga-Ra} which
reduces checking that $J/pJ$ is generated by the $p$-th
powers of its elements. However, the foregoing already
shows that $J/pB$ is generated by the $p$-th powers of its
elements; combining with theorem \ref{th_not-an-F_p}(ii),
the contention follows easily.
\end{pfclaim}

Taking into account remark \ref{rem_branch-ideal}(i),
we see that claim \ref{cl_two-squares} already implies
the assertion, in case $B_0$ is an $\F_p$-algebra.
Thus, we may assume that $B_0$ is not an $\F_p$-algebra,
in which case -- again by claim \ref{cl_two-squares} --
it suffices to show that $J\cap\fr B=\fr\cdot J$, or
equivalently, that :
$$
\Tor^B_1(B/\fr B,B/J)=0.
$$
Now, the chart $\beta$ yields a ring homomorphism
$C:=\Z[P^{(\infty)}]\to B$, and we have a base change
spectral sequence for $\Tor$-functors (see \cite[Th.5.6.6]{We}) :
$$
E^2_{pq}:=\Tor^B_p(\Tor^C_q(C/\fr C,B),B/J)\Rightarrow
\Tor^C_{p+q}(C/\fr C,B/J).
$$
\begin{claim}\label{cl_let-the-sunshine-in}
$\Tor^C_q(C/\fr C,B)=0$ for every $q>0$.
\end{claim}
\begin{pfclaim} Set $\fr^{(n)}:=\fr\cap P^{(n)}$ and
$C_n:=\Z[P^{(n)}]$ for every $n\in\N$. We have natural
isomorphisms :
$$
\Tor^C_q(C/\fr C,B)\isom\colim_{n\in\N}\Tor^{C_n}_q(C_n/\fr_nC,B_n)
\qquad
\text{for every $q,n\in\N$}
$$
so we are reduced to showing that $\Tor^{C_n}_q(C_n/\fr_nC,B_n)=0$
for every $q,n\in\N$ with $q>0$. The latter follows from
propositions \ref{prop_second-crit}(ii) and
\ref{prop_fla-criterion-point}(i).
\end{pfclaim}

From claim \ref{cl_let-the-sunshine-in} we see that $E^2_{pq}=0$
for every $q>0$ and every $p\in\N$. Thus, we deduce a natural
isomorphism :
$$
\Tor^B_p(B/\fr B,B/J)\isom\Tor^C_p(C/\fr C,B/J)
\qquad
\text{for every $p\in\N$}.
$$
By the same token, the projection $C\to \bar C:=C/pC$ induces
a base change spectral sequence :
$$
F^2_{ij}:=\Tor^{\bar C}_i(\Tor_j^C(C/\fr C,\bar C),B/J)\Rightarrow
\Tor_{i+j}^C(C/\fr C,B/J)
$$
(notice that $B/J$ is a $\bar C$-algebra, since $p\in J$) and
$\Tor_j^C(C/\fr C,\bar C)=0$ for every $j>0$, since $C/\fr C$
has no $p$-torsion. Hence $F^2_{ij}=0$ whenever $j>0$, and we
get an isomorphism :
$$
\Tor^{\bar C}_i(\bar C/\fr\bar C,B/J)\isom\Tor^C_i(C/\fr C,B/J)
\qquad
\text{for every $i\in\N$}.
$$
Notice now that $\bar C$ is a perfect $\F_p$-algebra, since
$P^{(\infty)}$ is $p$-divisbile, and with $I_{\bar C}:=\fr\bar C$
we have $\bar C/\fr\bar C=\bar C/I_{\bar C}^{\lceil 0\rceil}\bar C$.
Then $\Tor^{\bar C}_i(\bar C/\fr\bar C,B/J)=0$ for every $i>0$,
by proposition \ref{prop_perf-and-Tors}, whence the contention.

(ii): It is clear that $I_1\cap I_2$ is a branch ideal;
we have inclusions
$$
(I_1\cap I_2)^2\subset I_1I_2\subset I_1\cap I_2
$$
and from (i) we know that the first ideal in this chain
coincides with the third, so $I_1I_2=I_1\cap I_2$.

Next, for $k=1,2$ let us write $I_k=J_k\cap\br_k B$,
where $J_k\subset B$ and $\fr_k\subset P^{(\infty)}$ are
radical ideals and $p\in J_k$. Set $Z_k:=\Spec\,B/J_k$,
$Z'_k:=\Spec\,B/\fr_k B$ for $k=1,2$; clearly
$$
Z:=\Spec\,B/(I_1+I_2)=(Z_1\cup Z'_1)\cap(Z_2\cup Z'_2).
$$
Now let $J\subset B$ be the largest ideal such that
$\Spec\,B/J=(Z_1\cap Z_2)\cup(Z_1\cap Z'_2)\cup(Z'_1\cap Z_2)$.
It follows easily that
$$
Z=(\Spec\,B/J)\cup(\Spec\,B/(\fr_1\cup\fr_2)B).
$$ 
Now, $J$ and $\fr_1\cup\fr_2$ are radical ideals, and
$p\in J$, so we deduce that the radical of $I_1+I_2$
is a branch ideal, and it remains only to show that
$I_1+I_2$ is a radical ideal. However, notice the
inclusions
$$
(\fr_1 J_1+\fr_2 J_2)^2\subset
I:=(J_1+J_2)(J_1+\fr_2 B)(\fr_1B+J_2)(\fr_1\fr_2B)\subset
(\fr_1 J_1+\fr_2 J_2).
$$
Especially, $\fr_1 J_1+\fr_2 J_2$ is contained in the
radical $I'$ of $I$,  Since we already know that
$\fr_kJ_k=I_k$ for $k=1,2$, it suffices to check that
each of the four factors of the middle term in this
chain of inclusions is a branch ideal : indeed, in
this case the same holds for $I$, and therefore $I=I'$,
so finally $I=I_1+I_2$.

$\bullet$\ \
First, set $J:=J_1+J_2$, and let $\bar J\subset B/pB$
be the image of $J$; since $p\in J$, in order to see
that $J$ is a branch ideal, it suffices to check that
$J$ is a radical ideal; then it suffices to prove that
the same holds for $\bar J$, and the latter assertion
will follow, if we show that $\Phi_B^{-1}\bar J=\bar J$,
where $\Phi_B:B/pB\to B/pB$ is the Frobenius endomorphism.
Now, say that $x\in B/pB$ and $x^p\in\bar J$; then
$x^p=x_1+x_2$ for some $x_i\in \bar J_i:=J_i/pB$ ($i=1,2$),
and since $\Phi_B$ is surjective and $\bar J_i$ is a
radical ideal, we can write $x_i=y_i^p$ for some
$y_i\in\bar J_i$. Consequently $x^p=(y_1+y_2)^p$, so
$x-y_1-y_2$ is nilpotent; but clearly $\bar J$ contains
the nilradical of $B/pB$, whence the claim.

$\bullet$\ \
Next we consider $J':=J_1+\fr_2B$; as in the previous case,
since $p\in J$, we see that $J'$ is a branch ideal, provided
$J'/pB$ is a radical ideal of $B/pB$. However, say that
$x\in B/pB$ and $x^p=y+z$, for some $y\in\bar J_1$ and
$z\in\fr_2B/pB$; arguing as in the foregoing case, we
may write $y=u^p$ for some $u\in\bar J_1$. Likewise,
since both $\Phi_B$ and the Frobenius endomorphism of
$P^{(\infty)}$ are surjective, we see that $z=w^p$ for
some $w\in\fr_2B/pB$, so $x-u-w$ is nilpotent, and
especially, it lies in $\bar J_1$, whence the claim.
The same argument applies as well to the factor $\fr_1B+J_2$.

$\bullet$\ \
Lastly, since $\fr_1$ and $\fr_2$ are both radical
ideals of $P^{(\infty)}$, the same holds for $\fr_1\fr_2$,
so $\fr_1\fr_2 B$ is a branch ideal.
\end{proof}

Proposition \ref{prop_branch-setup}(i) says that
the pair $(B,I)$ is a basic setup, in the sense of
\cite[(2.1.1)]{Ga-Ra}, so we may consider the associated
categories of almost modules and almost algebras.

\sset\subsubsection{}\label{subsec_strata-tower}
Henceforth we shall restrict to the case where $B_0$ is
not an $\F_p$-algebra, so $B_0$ is a local integral
domain whose field of fractions has characteristic zero,
and whose residue field has characteristic $p$. For every
$n\in\N$, set
$$
(X_n)_{/p}:=\Spec\,B_n/pB_n
\qquad
V_n:=X_n\setminus(X_n)_{/p}
\qquad
V_{n,\tr}:=(X_n,\underline M)_\tr\cap V_n.
$$
We consider now a finite morphism $\phi_0:Y_0\to X_0$ with $Y_0$
a normal scheme, such that $\phi_0$ maps each connected component
of $Y_0$ onto $X_0$, and such that the restriction
$\phi^{-1}_0V_{0,\tr}\to V_{0,\tr}$ is a finite \'etale covering.
For every $n\in\N$, we let $Y_n$ be the normalization of $X_n$
in $\phi^{-1}V_{0,\tr}\times_{X_0}X_n$, and denote by
$\phi_n:Y_n\to X_n$ the resulting finite morphism
(\cite[Lemma 1, p.262]{Mat}), by
$$
\phi:Y\to X
$$
the limit of the system $(\phi_n~|~n\in\N)$, and by
$W_n\subset X_n$ the {\em \'etale locus of $\phi_n$},
{\em i.e.} the largest open subset such that the restriction
$\phi_n^{-1}W_n\to W_n$ of $\phi_n$ is an \'etale covering
(lemma \ref{lem_replace}(iii) and claim \ref{cl_trivial-Zorni}).
Let also $(X_{n,\fp}~|~\fp\in\Spec\,P^{(n)})$ be the
logarithmic stratification of $(X_n,\underline M{}_n)$ (see
\eqref{subsec_log-stratif}) and set
$V_{n,\fp}:=X_{n,\fp}\cap V_n$ for every $n\in\N$ and every
$\fp\in\Spec\,P^{(n)}$. We call
$(V_{n,\fp}~|~\fp\in\Spec\,P^{(n)})$ the {\em logarithmic
stratification\/} of $V_n$.

\begin{lemma}\label{lem_log-stratif}
With the notation of \eqref{subsec_strata-tower}, the
following holds :
\begin{enumerate}
\item
$V_n\setminus W_n$ is a union of strata of the logarithmic
stratification of $V_n$, for every $n\in\N$.
\item
$X_n\times_{X_m}W_m=W_n$ for every sufficiently large
$m\in\N$ and every $n\geq m$.
\item
For $m$ as in {\em (ii)}, set $W:=X\times_{X_m}W_m$, and
let $I\subset B$ be the largest ideal such that
$\Spec\,B/I=X\!\setminus\!W$. Then $I$ is a branch
ideal of $B$.
\end{enumerate}
\end{lemma}
\begin{proof}(i): Since every $X_{n,\fp}$ is irreducible
(corollary \ref{cor_logar-strata}(iii)), the same holds for
each $V_{n,\fp}$. Hence, suppose that $V_n\setminus W_n$ contains
a point $z$ of $V_{n,\fp}$, for some $\fp\in\Spec\,P^{(n)}$; we
have to show that in this case $V_{n,\fp}\cap W_n=\emptyset$,
and since $W_n$ is an open subset of $X_n$, it suffices to
check that the generic point $\eta$ of $V_{n,\fp}$ does not
lie in $W_n$.
However, let $\bar\eta$ be any geometric point localized
at $\eta$; then $\eta$ lies in $W_n$ if and only if
the induced morphism $\phi_n\times_{X_n}X_n(\bar\eta)$
is \'etale (claim \ref{cl_include-pts}). Let also $\bar z$
be a geometric point localized at $z$, and denote by
$\underline M{}_{n,\bar z}$ (resp. $\underline M{}_{n,\bar\eta}$)
the stalk at $\bar z$ (resp. $\bar\eta$) of the logarithmic
structure of $(X_n(\bar z),\underline M{}_n(\bar z))$ (resp.
of $(X_n(\bar\eta),\underline M{}_n(\bar\eta))$). Any choice
of a strict specialization morphism
$s:X_n(\bar\eta)\to X_n(\bar z)$ induces a strict specialization
map
$\bar\sigma:\underline M{}_{n,\bar z}\to\underline M{}_{n,\bar\eta}$
(see \eqref{subsec_strict-special}) that extends the
specialization map
$\sigma:\underline M{}_{n,z}\to\underline M{}_{n,\eta}$.
Since the natural maps
$\underline  M{}_{n,z}^\sharp\to\underline  M{}_{n,\bar z}^\sharp$
and
$\underline  M{}_{n,\eta}^\sharp\to\underline  M{}_{n,\bar\eta}^\sharp$
are isomorphisms (see \eqref{eq_general-iso}), and
$\sigma^{\sharp\gp}$ is an isomorphism (since $z$ and $\eta$
lie in the same stratum), the same holds for
$\bar\sigma{}^{\sharp\gp}$, and therefore the induced
map $\underline M{}^{\gp\vee}_{n,\bar\eta}\to
\underline M{}^{\gp\vee}_{n,\bar z}$ is bijective.
Pick any geometric point $\xi$ of $X_n(\bar z)$ localized
at the maximal point, and lift $\xi$ to a geometric point
$\xi_\eta$ of $X_n(\bar\eta)$; in light of
\eqref{subsec_punctured}, we conclude that $s$ induces
an isomorphism
$$
\pi_1((X_n(\bar\eta),\underline M{}_n(\bar\eta))_\et,\xi_\eta)
\isom\pi_1((X_n(\bar z),\underline M{}_n(\bar z))_\et,\xi).
$$
Therefore, $\phi_n\times_{X_n}X_n(\bar\eta)$ is \'etale
if and only if it is a trival covering, if and only if
the same holds for $\phi_n\times_{X_n}X_n(\bar z)$,
if and only if $z\in W_n$, whence the assertion. 

(ii): For every $n\in\N$, set
\set\begin{equation}\label{eq_Z_n-for-later}
\fZ_n:=\{\fp\in\Spec\,P^{(n)}~|~V_{n,\fp}\neq\emptyset
\quad\text{and}\quad V_{n,\fp}\cap W_n=\emptyset\}.
\end{equation}
The continuous map $\omega_n:\Spec\,P^{(n+1)}\to\Spec\,P^{(n)}$
induced by the inclusion $P^{(n)}\to P^{(n+1)}$ sends $\fZ_{n+1}$
into $\fZ_n$, for every $n\in\N$. On the other hand, $\omega_n$
is also a bijection of finite sets (lemmata \ref{lem_face}(iii)
and \ref{lem_Kummer-fans}(i)); we conclude that $\omega_n$
restricts to a bijection $\fZ_{n+1}\isom\fZ_n$ for every
sufficiently large integer $n\in\N$. For every $n,m\in\N$
with $n\geq m$, let $g_{n,m}:X_n\to X_m$ be the transition
morphism in the maximal tower \eqref{eq_max-tower}, in light
of (i), it follows that
$$
g^{-1}_{n,m}(V_m\cap W_m)=V_n\cap W_n
\qquad
\text{for every sufficiently large $m\in\N$ and every $n\geq m$}.
$$
On the other hand, the restriction $(X_n)_{/p}\to(X_m)_{/p}$
of $g_{n,m}$ is radicial and surjective for every such
$n,m\in\N$, so the underlying continuous map is a homeomorphism;
since the underlying topological spaces are noetherian, we see
as well that
$$
g^{-1}_{n,m}((X_m)_{/p}\cap W_m)=(X_n)_{/p}\cap W_n
\qquad
\text{for every sufficiently large $m\in\N$ and every $n\geq m$}.
$$
Summing up, the assertion follows.

(iii): The assertion follows easily from (i) and
\eqref{eq_radical-identity}.
\end{proof}

\sset\subsubsection{}\label{subsec_keep-this-on}
Let $I\subset B$ be a fixed branch ideal, and consider a
pair $(X,Z)$ with $Z\subset\Spec\,B/I$ and such that $Z$
is constructible in $X$. Clearly such a pair is normal
(see definition \ref{def_al-pure}(ii)), and we aim
to show that $(X,Z)$ is almost pure for the almost
structure given by the basic setup $(B,I)$ supplied
by proposition \ref{prop_branch-setup}. This will be
achieved in several steps.

To begin with, set $U:=X\setminus Z$, and let $\cA$ be
any \'etale almost finitely presented $\cO_{\!U}^a$-algebra.
Set also $U_{\!I}:=X\setminus\Spec\,B/I$; then $\cA_{|U_{\!I}}$
is a finite \'etale $\cO_{U_{\!I}}$-algebra,
$U_{\!I}=X\times_{X_0}U_{\!I,0}$ for some open subset
$U_{\!I,0}\subset X_0$ (remark \ref{rem_branch-ideal}(iv)),
and a simple inspection shows that
$$
V_{0,\tr}\subset U_{\!I,0}.
$$
Let $\psi:U_{\!I}\to U_{\!I,0}$ be the natural projection;
by \cite[Ch.IV, Prop.17.7.8(ii)]{EGA4} and
\cite[Ch.IV, Th.8.8.2(ii)]{EGAIV-3}, and by virtue
of remark \ref{rem_shift-sequence}, we may assume --
after replacing $(X_0,\underline M{}_0)$ by
$(X_n,\underline M{}_n)$ for some sufficiently large
$n\in\N$ -- that there exists a coherent \'etale
$\cO_{\!U_{\!I,0}}$-algebra $\cA_0$ with an isomorphism
$\psi^*\cA_0\isom\cA_{|U_I}$ of $\cO_{\!U_I}$-algebras.
Denote by $Y_0$ the normalization of $X_0$ in
$\Spec\,\cA_0(U_{\!I,0})$; the resulting morphism
$\phi_0:Y_0\to X_0$ is finite (lemma \ref{lem_Nagatas}(i))
and induces an isomorphism
$(\phi_{0*}\cO_{Y_0})_{|U_{\!I,0}}\isom\cA_0$. Especially,
$\phi_0$ is a morphism of the type contemplated in
\eqref{subsec_strata-tower}, and the resulting morphism
$\phi:Y\to X$ induces an isomorphism
$(\phi_*\cO_Y)_{|U_{\!I}}\isom\cA_{|U_{\!I}}$. However, we
have as well $\cA=(\cA_{|U_{\!I}})^\nu$ (lemma
\ref{lem_pure-almost-crit}(i)), and since $Y$ is normal,
there follows an isomorphism of $\cO^a_{\!U}$-algebras
$(\phi_*\cO^a_Y)_{|U}\isom\cA$. Then proposition
\ref{prop_almost-pure-crit} and \cite[Lemma 8.2.28]{Ga-Ra}
show that $(X,Z)$ is almost pure if and only if the
resulting $\cO^a_{\!X}$-algebra $\phi_*\cO^a_Y$ is weakly
unramified, for every such $\phi_0$.

\sset\subsubsection{}\label{subsec_was-a-dot}
Keep the notation of \eqref{subsec_keep-this-on},
and {\em suppose now that $B_0$ is strictly henselian};
especially, the $N$-torsion subgroup $\mu_N$ of
$B^\times_0$ has cardinality equal to $N$, for every
$N>0$ such that $(N,p)=1$.
Let $k>0$ be any integer and write $k=p^s\cdot q$, with
$s,q\in\N$ and $(q,p)=1$; set $Q:=P$, let $\nu:P\to Q$
be the $k$-Frobenius map, define
$$
C'_0:=Q\otimes_PB_0
\qquad
C_0:=C'_0\otimes_{B_0}B''_s
$$
and endow $X'_0:=\Spec\,C_0$ with the log structure
$\underline M{}'_0$ deduced from the natural map $Q\to C_0$.
Notice that $C_0$ is strictly henselian and the Frobenius
endomorphism of $C_0/pC_0$ is still a finite map
(claim \ref{cl_stability}(i)); also, since the induced
map $A\to C'_0/\fm_QC'_0$ is an isomorphism, the image
of the sequence $(f_1,\dots,f_r)$ in $C_0$ is still maximal,
in the sense of \eqref{subsec_more-notation}, and therefore
$C_0/\fm_QC_0$ is still a regular local ring, by corollary
\ref{cor_regu-criterion}. Arguing as in
\eqref{subsec_another-Frobenius}, we deduce that
$(X'_0,\underline M{}'_0)$ is a regular log scheme, and
we may consider the maximal tower
$$
((X'_n,\underline M{}'_n)~|~n\in\N)
$$
associated with the chart $Q\to C_0$ and the maximal
sequence $(f^{1/p^s}_1,\dots,f^{1/p^s}_r)$ (see remark
\ref{rem_afterthought}(iii)). So, $X'_n=\Spec\,C_n$
for a finite $C_0$-algebra $C_n$, and $\underline M{}'_n$
is given by a chart $Q^{(n)}\to C_n$, where $Q^{(n)}$ is
a submonoid of $Q_\Q$ containing $Q$, for every $n\in\N$.
For every $n\in\N$, we set
$V'_{n,\tr}:=V_{0,\tr}\times_{X_0}X'_n$, we let $Y'_n$ be
the normalization of $X'_n$ in $Y_0\times_{X_0}V'_{n,\tr}$,
we denote by
$$
\phi'_n:Y'_n\to X'_n
\qquad\text{and}\qquad
h_n:(X'_n,\underline M{}'_n)\to(X_n,\underline M{}_n)
$$
the resulting finite morphisms, and by
$\phi':Y'\to X'$ (resp. $h:(X',\underline M')\to(X,\underline M)$)
the limit of the system of morphisms $(\phi'_n~|~n\in\N)$ (resp.
the limit of the system $(h_n~|~n\in\N)$).
Recall also that the induced morphism $V'_{0,\tr}\to V_{0,\tr}$
is a torsor for the finite abelian group
$$
G:=\Hom_\Z(Q^\gp,\mu_k)
$$
(see \eqref{subsec_Second-pairing}). Set
$$
C:=\bigcup_{n\in\N}C_n
\qquad
Q^{(\infty)}:=\bigcup_{n\in\N}Q^{(n)}.
$$
With this notation, we have
$$
Q^{(\infty)}=Q\otimes_PP^{(\infty)}
\qquad\text{and}\qquad
C=Q^{(\infty)}\otimes_PB''.
$$
Let also write $Y=\Spec\,D$ (resp. $Y'=\Spec\,D'$) for a
$B$-algebra $D$ (resp. for a $C$-algebra $D'$).
The morphism $h:(X',\underline M')\to(X,\underline M)$
shall be called the {\em standard covering} of
$(X,\underline M)$ of degree $k$.

\begin{lemma}\label{lem_reduce-to-p-branch}
In the situation of \eqref{subsec_was-a-dot}, set
$U':=X'\setminus h^{-1}Z$. The induced diagram
$$
\xymatrix{
\cO^a_{\!X}\Et_\mathrm{fr} \ar[r]^-\rho \ar[d] &
\cO^a_{\!U}\Et_\mathrm{fr} \ar[d] \\
\cO^a_{\!X'}\Et_\mathrm{fr} \ar[r]^-{\rho'} &
\cO^a_{\!U'}\Et_\mathrm{fr}
}$$
is $2$-cartesian (notation of definition {\em\ref{eq_resr_etale}}).
\end{lemma}
\begin{proof} A simple inspection shows that $h_0$ factors through
a morphism $(X'_0,\underline M{}'_0)\to(X_s,\underline M{}_s)$, and
indeed we have natural isomorphisms of $(X',\underline M')$-schemes
$$
(X'_n,\underline M{}'_n)\isom(X'_0,\underline M{}'_0)
\times_{(X_s,\underline M{}_s)}(X_{n+s},\underline M_{n+s})
\qquad
\text{for every $n\in\N$}
$$
that identify the transition morphisms
$(X'_{n+1},\underline M{}'_{n+1})\to(X'_n,\underline M{}'_n)$
with the base change of the corresponding morphisms for the
tower $((X_{n+s},\underline M{}_{n+s})~|~n\in\N)$. In view
of remark \ref{rem_shift-sequence}, we may then replace
$(X_0,\underline M_0)$ by $(X_s,\underline M{}_s)$, and $P$ by
$P^{(s)}$, and assume that $(k,p)=1$, so also $(o(G),p)=1$.

Suppose now that $(\phi'_*\cO^a_{Y'})_{|U'}$ is in the
essential image of the restriction functor $\rho'$; by
proposition \ref{prop_almost-pure-crit}, this means that
$D^{\prime a}$ is an \'etale almost finitely presented
$(C,IC)^a$-algebra. Then, taking into account the discussion of
\eqref{subsec_keep-this-on}, and recalling that $\rho$ and
$\rho'$ are fully faithful (lemma \ref{lem_pure-almost-crit}(ii)),
we are reduced to checking that $D^a$ is an \'etale almost
finitely presented $(B,I)^a$-algebra.

However, the action of $G$ on $V'_{0,\tr}$ is inherited by $C$ and
$D'$, and clearly $C^G=B$ and $D^{\prime G}=D$; by corollary
\ref{cor_horiz}(ii), it suffices therefore to show that the
action of $G$ on $D^{\prime a}$ is horizontal. Denote
$$
G\to\Aut(C)
\quad :\quad
\chi\mapsto\rho_\chi
$$
the action of $G$, let $\chi\in G$ be any element, and
$J_\chi\subset C$ the ideal generated by all elements
of the form $\rho_\chi(c)-c$, for $c$ ranging over all
elements of $C$. By the discussion in
\eqref{subsec_Second-pairing} we get
$$
\rho_\chi((q\otimes y)\otimes b)=\chi(q)\cdot(q\otimes y)\otimes b
\qquad
\text{for every $q\in Q$, $y\in P^{(\infty)}$ and $b\in B''$}.
$$
Hence, $J_\chi$ is the ideal generated by all elements of
the form $(1-\chi(q))\cdot q\otimes y\otimes 1$, for all
$q\otimes y\in Q^{(\infty)}$. However, since $\chi(q)\in\mu_k$
and $(k,p)=1$, it is easily seen that $1-\chi(q)$ either
vanishes, or else it is invertible in $B_0$. Thus, denote
by $\fq_\chi\subset Q^{(\infty)}$ the ideal generated by all
elements of the form $q\otimes y$, with $\chi(q)\neq 1$;
it follows that
$$
J_\chi=\fq_\chi\cdot C.
$$

\begin{claim}\label{cl_radical-ideals}
$\fq_\chi$ is a radical ideal.
\end{claim}
\begin{pfclaim} Indeed, say that $(q\otimes y)^n\in\fq_\chi$,
so there exist elements $q_1\otimes y_1$ and $x_1$ of
$Q^{(\infty)}$ such that
$(q\otimes y)^n=(q_1\otimes y_1)\cdot x_1$ and $\chi(q_1)\neq 1$.
We may assume that $n=p^t$ for some integer $t\in\N$, and since
$Q^{(\infty)}$ is uniquely $p$-divisible, we may write
$q_1\otimes y_1=(q_2\otimes y_2)^n$ and $x_1=x_2^n$ for
some elements $q_2\otimes y_2$ and $x_2$ of $Q^{(\infty)}$,
and then $q\otimes y=(q_2\otimes y_2)\cdot x_2$; however,
clearly $\chi(q_2)\neq 1$, whence the claim.
\end{pfclaim}

In view of claim \ref{cl_radical-ideals}, we may write
$\fq_\chi=\fq_1\cap\cdots\cap\fq_s$, for certain
$\fq_1,\dots,\fq_s\in\Spec\,Q^{(\infty)}$ (lemma \ref{lem_radical}).
Set $\fQ:=\{\fq_1,\dots,\fq_s\}$; we also know that
$(\fq\cap Q^{(n)})\cdot C_n$ is a prime ideal of $C_n$,
for every $\fq\in\fQ$ and every $n\in\N$
(corollary \ref{cor_same-height}(i)), therefore $\fq C$
is a prime ideal of $C$, for every $\fq\in\fQ$.
Especially, the induced map
\set\begin{equation}\label{eq_use-the-tag}
C/J_\chi\to\prod_{\fq\in\fQ}C/\fq C
\end{equation}
is injective. Since, by assumption, $D^{\prime a}$ is a flat
$C^a$-algebra, it follows that the map of almost modules
$D^{\prime a}\otimes_C\eqref{eq_use-the-tag}$ is a
monomorphism. Summing up, we are reduced to checking that
$\chi$ acts trivially on $(D'/\fq D')^a$, for every
$\fq\in\fQ$. However, set
$$
X'_\fq:=\Spec\,C/\fq C
\quad\text{and}\quad
U'_\fq:=X'_\fq\times_XU_I
\qquad
\text{for every $\fq\in\fQ$}.
$$
Suppose first that $U'_\fq=\emptyset$; in that case,
set $\fp:=\fq\cap P^{(\infty)}$, and notice that $\fq$
is the radical of the ideal $\fp\cdot Q^{(\infty)}$, so
$\fq C$ is the radical of $\fp C$, and
therefore $\Spec\,(C/\fp C)\times_XU_I=\emptyset$.
However, the induced morphism
$\Spec\,C/\fp C\to\Spec\,B/\fp B$ is surjective, so
$\Spec\,(B/\fp B)\cap U_I=\emptyset$; since
$\fp B$ is a prime ideal of $B$, the latter means
that $I\subset\fp B$, whence $(B/\fp B)^a=0$,
so $(D'/\fq D')^a$ vanishes as well, and the assertion
is trivial. If $U'_\fq\neq\emptyset$, set
$Y'_\fq:=Y'\times_{X'}X'_\fq$, let $\phi'_\fq:Y'_\fq\to X'_\fq$
be the induced morphism, and define
$$
\cD:=\phi_*\cO_Y
\qquad
\cD':=\phi'_*\cO_{Y'}
\qquad
\cD'_\fq:=\phi'_{\fq*}\cO_{Y'_\fq}.
$$
On the one hand, by assumption $\cD_\fq^{\prime a}$ is a flat
$\cO^a_{\!X'_\fq}$-algebra; on the other hand, $X'_\fq$ is
reduced and irreducible, hence the restriction map
$\Gamma(X'_\fq,\cO_{\!X'_\fq})\to\Gamma(U'_\fq,\cO_{\!X'_\fq})$
is injective. Consequently, the restriction map
$$
(D'/\fq D')^a=\Gamma(X'_\fq,\cD^{\prime a}_\fq)
\to\Gamma(U'_\fq,\cD^{\prime a}_\fq)
$$
is a monomorphism; so, we are reduced to checking that
$\chi$ acts trivially on $\Gamma(U'_\fq,\cD^{\prime a}_\fq)$.
To this aim, we remark that the open subset
$U'_{\!I}:=X'\times_XU_{\!I}$ is stable under the action of
$G$, and the restriction $\cD'_{|U'_{\!I}}$ of $\cD'$ to the
open subset $U'_{\!I}$ is isomorphic to $(h^*\cD)_{|U'_{\!I}}$
(\cite[Ch.IV, Prop.17.5.8(iii)]{EGA4}), so the action of
$G$ on $\cD'_{|U'_{\!I}}$ is horizontal, and the assertion
follows easily. 
\end{proof}

\sset\subsubsection{}\label{subsec_kezp}
In the situation of \eqref{subsec_was-a-dot}, define
$\fZ_n$ as in \eqref{eq_Z_n-for-later}, for every
$n\in\N$. As already remarked, in view of lemma
\ref{lem_log-stratif}(i,ii) we may assume that the
map $\Spec\,P^{(n+1)}\to\Spec\,P^{(n)}$ sends $\fZ_{n+1}$
bijectively onto $\fZ_n$, for every $n\in\N$. Set
$V'_0:=V_0\times_{X_0}X'_0$, and let
$$
(V'_{0,\fp}~|~\fp\in\Spec\,Q)
\qquad\text{and}\qquad
(V_{0,\fp}~|~\fp\in\Spec\,P)
$$
be the logarithmic stratifications of $V'_0$ and respectively
$V_0$, defined as in \eqref{subsec_strata-tower}; notice that
the map $\nu^*:\Spec\,Q\to\Spec\,P$ induced by $\nu$ is
bijective (lemma \ref{lem_Kummer-fans}(i)), and clearly
$$
V'_{0,\fp}=h_0^{-1}V_{0,\nu^*(\fp)}
\qquad
\text{for every $\fp\in\Spec\,Q$}.
$$
Let $\fZ'_0:=\nu^{*-1}\fZ_0$, and for every $\fp\in\fZ'_0$, let
$\eta'_\fp$ (resp. $\eta_\fp$) denote the generic point
of $V'_{0,\fp}$ (resp. of $V_{0,\nu^*(\fp)}$), pick geometric
points $\bar\eta{}'_\fp$ and $\xi'_\fp$ localized
respectively at $\eta'_\fp$ and at a point of
$V'_{0,\tr}(\bar\eta{}'_\fp):=V'_{0,\tr}\times_{X'_0}X'_0(\bar\eta{}'_\fp)$.
Denote by $\bar\eta_\fp$ the image of $\bar\eta{}'_\fp$
in $V_{0,\nu^*(\fp)}$, and by $\xi_\fp$ the image of $\xi'_\fp$
in $V_{0,\tr}(\bar\eta_\fp):=V_{0,\tr}\times_{X_0}X_0(\bar\eta_\fp)$.
According to \eqref{subsec_funct-standard-shit} there follows,
for every integer $N>0$, a commutative diagram of groups
$$
\xymatrix{
\pi_1(V'_{0,\tr}(\bar\eta{}'_\fp)_\et,\xi'_\fp) \ar[r] \ar[d] &
\pi_1(V_{0,\tr}(\bar\eta_\fp)_\et,\xi_\fp) \ar[d] \\
\underline M{}_{0,\eta'_\fp}^{\prime\gp\vee}\otimes_\Z\mu_N \ar[r] &
\underline M{}_{0,\eta_\fp}^{\gp\vee}\otimes_\Z\mu_N
}$$
(where $\mu_N$ is the $N$-torsion subgroup of
$\kappa(\xi_\fp)^\times$) whose top arrow is induced by the
natural morphism
$V'_{0,\tr}(\bar\eta{}'_\fp)\to V_{0,\tr}(\bar\eta_\fp)$,
and whose bottom arrow is induced by
$\nu^{\gp\vee}:Q^{\gp\vee}\to P^{\gp\vee}$, {\em i.e.} by the
$k$-Frobenius map of $P^{\gp\vee}$, for every $\fp\in\fZ'_0$.
Now, the restriction
$$
\phi_\fp:Y_0(\bar\eta_\fp):=Y_0\times_{X_0}X_0(\bar\eta_\fp)
\to X_0(\bar\eta_\fp)
$$
of $\phi$ is a tamely ramified covering, hence the
action of $\pi_1(V_{0,\tr}(\bar\eta_\fp)_\et,\xi_\fp)$ on
$F_\fp:=\phi^{-1}_\fp(\xi_\fp)$ factors through a
group homomorphism
$$
\rho_\fp:
\underline M{}_{0,\eta_\fp}^{\gp\vee}\otimes_\Z\mu_N\to\Aut(F_\fp)
$$
for some sufficiently large $N\in\N$ (theorem
\ref{th_Abhyankar_log}).
We may then find $k\in\N$ such that the image of
$k\cdot P^{\gp\vee}$ in $\underline M{}_{0,\eta_\fp}^{\gp\vee}$
lies in the kernel of $\rho_\fp$, for every $\fp\in\fZ'_0$.
Especially, for this choice of $k$, the image of
$\underline M{}_{0,\eta'_\fp}^{\prime\gp\vee}\otimes_\Z\mu_N$
acts trivially on $\Aut(F_\fp)$ via $\rho_\fp$, for every
such $\fp$. Consequently,
$\pi_1(V'_{0,\tr}(\bar\eta{}'_\fp)_\et,\xi'_\fp)$
acts trivially on the fibres $\phi_0^{\prime-1}(\xi'_\fp)$
(by virtue of \eqref{eq_restrict-action}); after applying
lemma \ref{lem_log-stratif}(i) to the morphism $\phi'_0$,
we conclude that the \'etale locus of $\phi'_0$ contains
the whole of $V'_0$.

\begin{theorem}\label{th_log-regular-purity}
In the situation of \eqref{subsec_keep-this-on},
the pair $(X,Z)$ is almost pure, relative to the
basic setup $(B,I)$.
\end{theorem}
\begin{proof} Fix a geometric point $\bar x$ of $X$ localized
at the closed point $x$ of $X$, and let $B^\sh$ denote the
strict henselization of $B$ at $x$. In view of proposition
\ref{prop_pro-smooth-desc-pure}, in order to prove that
$(X,Z)$ is almost pure, it suffices to show that the pair
$(X(\bar x),Z(\bar x))$ is almost pure relative to the basic
setup $(B,I)$, or -- equivalently -- relative to the basic
setup $(B^\sh,IB^\sh)$. However, it is easily seen that
$IB^\sh$ is also a branch ideal of $B^\sh$, if the
latter is endowed with the chart $P^{(\infty)}\to B^\sh$
deduced from the given chart of $\underline M$.
Taking into account remark \ref{rem_reduce-to-sh}(i),
we may then replace $B$ by $B^\sh$, and assume
that $X_0$ and $X$ are strictly local.

Now, consider a finite morphism $\phi_0:Y_0\to X_0$
and the resulting morphism $\phi:Y\to X$ as in
\eqref{subsec_strata-tower}. The discussion of
\eqref{subsec_keep-this-on} shows that $(X,Z)$ is almost
pure if and only if $\phi_*\cO^a_Y$ is an \'etale
$\cO^a_{\!X}$-algebra, for every such $\phi_0$. Set
$V:=\Spec\,B[1/p]$; in view of \eqref{subsec_kezp} and
lemma \ref{lem_reduce-to-p-branch}, we may then assume
that $V_0\subset U_0$, and therefore $(\phi_*\cO_Y)_{|V}$
is a finite \'etale $\cO_{\!X}$-algebra. On the other
hand, let $I'\subset B$ be the radical of the ideal
$I+pB$, and set $Z':=Z\!\setminus\!V$; clearly, $I'$ is a
branch ideal. It then suffices to show that $\phi_*\cO^a_Y$
is an \'etale $\cO^a_{\!X}$-algebra, for the almost structure
given by the new setup $(B,I')$. Thus, we may assume
that $p\in I$, in which case the assertion follows from
theorems \ref{th_not-an-F_p}(iii) and
\ref{th-alm-purity-form-perfectoid}.
\end{proof}

\begin{theorem}\label{th_B-ind-meas-reglog}
In the situation of \eqref{subsec_another-Frobenius}, the
ring $B$ is ind-measurable.
\end{theorem}
\begin{proof} (See remark \ref{rem_ind-measure}(i) for the
definition of ind-measurable ring.) One argues as in the proof
of theorem \ref{th_B-sh-is-measurable}, with some simplifications.
We have to exhibit a sequence $(d_n~|~n\in\N)$ of normalizing
factors fulfilling conditions (a) and (b) of definition
\ref{def_ind-measur}, where the $\lambda_n$ occurring in
{\em loc.cit.} is meant to be the usual length function for
finitely generated $B_n$-modules supported at the
closed point $x_n$ of $X_n$. Now, fix $n\in\N$, set
$$
\bT_\R:=P^\gp_\R/P^{(n)\gp}
$$
and endow $\bT_\R$ with its invariant measure $d\mu_n$ of
total volume equal to $1$. For every $\gamma\in P^\gp_\R$,
let $[\gamma]\in \bT_\R$ be the equivalence class of
$\gamma$; notice that the $P^{(n)}$-module
\set\begin{equation}\label{eq_S-is-for-module}
S_{[\gamma]}:=\gamma P^{(n)\gp}\cap P_\Q
\end{equation}
is finitely generated (proposition \ref{prop_Gordon}(ii))
and depends only on the class $[\gamma]$, and for any given
finitely generated $B_n$-module $M$ supported at $x_n$,
consider the function
$$
l_M:\bT_\R\to\N
\qquad
[\gamma]\mapsto\lambda_n(S_{[\gamma]}\otimes_{P^{(n)}}M).
$$
Let $e_1,\dots,e_r$ be a basis of the free
$\Z$-module $P^\gp$, and define $\Omega_n\subset P^\gp_\R$
as in the proof of theorem \ref{th_B-sh-is-measurable}, so
that $\Omega_n$ is a fundamental domain for the lattice
$P^{(n)\gp}$, and $0$ lies in the interior of $\Omega_n$.
Denote also by $\Sigma\subset\bT_\R$ the image of
$P_\R\cap\Omega_n$.

\begin{claim}\label{cl_partition}
There is a partition of $\bT_\R$ into finitely many measurable
subsets $\Theta_1,\dots,\Theta_t$, independent of $M$, such
that :
\begin{enumerate}
\item
For every $\gamma,\lambda\in\Omega_n$, the classes
$[\gamma],[\lambda]\in\bT_\R$ lie in the same $\Theta_i$
if and only if $\gamma^{-1}S_{[\gamma]}=\lambda^{-1}S_{[\lambda]}$.
\item
Especially, $l_M$ restricts to a constant function
on each $\Theta_i$.
\item
Let $\Theta\in\{\Theta_1,\dots,\Theta_t\}$ be the subset
containing $[0]\in\bT_\R$; then $\Theta\cap\Sigma$ has
measure $>0$.
\end{enumerate}
\end{claim}
\begin{pfclaim} According to proposition
\ref{prop_linear-constr-part}(i,iii), the set
$\cS:=\{\gamma^{-1}S_{[\gamma]}~|~\gamma\in\Omega_n\}$
is finite, and for every non-empty $S\in\cS$, the set
$\{\gamma\in\Omega_n~|~\gamma^{-1}S_{[\gamma]}=S\}$ is the
intersection of $\Omega_n$ with a $\Q$-linearly constructible
subset. It follows that the same must hold also in case
$S=\emptyset$. The image in $\bT_\R$ of any such $\Q$-linearly
constructible subset is obviously measurable, whence (i).
Moreover, in view of our choice of $\Omega$, assertion (ii)
follows easily from claim \ref{cl_C_eps-lies}. 
\end{pfclaim}

Let $m\geq n$ be any integer; since $B''_m$ is a free $B''_n$-module
of rank $p^{r(m-n)}$ (notation of \eqref{subsec_more-notation}), we may
compute :
$$
\lambda_m(B_m\otimes_{B_n}M)=
\frac{p^{r(m-n)}}{[\kappa(x_m):\kappa(x_n)]}\cdot
\sum_{[\gamma]\in P^{(m)\gp}/P^{(n)\gp}}l_M([\gamma]).
$$
However, lemma \ref{lem_Frob-onto} easily implies that
$\kappa(x_{n+1})^p=\kappa(x_n)$ for every $n\in\N$, whence
$$
[\kappa(x_{n+1}):\kappa(x_n)]=p^{e_n}
\qquad
\text{where $e_n:=\Omega^1_{\kappa(x_n)/\Z}$}
$$
by virtue of \cite[Ch.IV, Th.21.4.5]{EGAIV}. But by the same
token, the field $\kappa(x_m)$ is isomorphic to $\kappa(x_n)$
for every $m\geq n$, so $e_n$ is actually independent of $n$,
and we get
$$
[\kappa(x_m):\kappa(x_n)]=p^{e_0(m-n)}
\qquad
\text{for every $m\geq n$}.
$$
Therefore, set
$$
d_n:=p^{n\cdot\dim B_0}
\qquad
\text{for every $n\in\N$}.
$$
We claim that $(d_n~|~n\in\N)$ is a suitable sequence of
normalizing factors for $B$. Indeed, claim \ref{cl_partition}(ii)
says that $l_M$ is a measurable function on $\bT_\R$, and the
foregoing, together with the discussion of \eqref{subsec_local-max}
implies that :
$$
\lambda(B\otimes_{B_n}M):=
\lim_{m\to+\infty}d_m^{-1}\cdot\lambda_m(B_m\otimes_{B_n}M)=
d^{-1}_n\int_{\bT_\R}l_Md\mu_n
$$
(recall that $\dim P=\rk_\Z P^\gp$, by corollary
\ref{cor_consequent}(i)), so condition (a) holds for this
choice of factors. Next, fix $\eps>0$, let $N\to N'$ be a
surjection of finitely generated $B_n$-modules supported
at $x_n$, and suppose that
$d_n^{-1}(\lambda_n(N)-\lambda_n(N'))\geq\eps$. Since
$\lambda_n(N)=l_N(0)$ (and likewise for $N'$), we deduce that
$$
\lambda(B\otimes_{B_n}N)-\lambda(B\otimes_{B_n}N')
\geq\eps\cdot\int_\Theta d\mu_n
$$
where $\Theta\in\{\Theta_1,\dots,\Theta_t\}$ is the unique
subset of $\bT_\R$ such that $[0]\in\Theta$, so the volume of
$\Theta$ is $>0$, by claim \ref{cl_partition}(iii). This
shows that condition (b) holds as well, and concludes the
proof of the theorem.
\end{proof}

\subsection{The direct summand conjecture}
\label{sec_dir-summand-conj}
In the following paragraphs we prove a generalization of
the direct summand conjecture for finite injective extensions
of log-regular rings.

To begin with, let $A$ be a ring, and consider a short exact
sequence of $A$-modules
$$
\Sigma\qquad : \qquad
0\to M'\xrightarrow{f}M\xrightarrow{g}M''\to 0.
\qquad\qquad\qquad
$$
We say that $\Sigma$ is {\em universally exact} if the sequence
$X\otimes_A\Sigma$ is short exact for every $A$-module $X$.

\begin{lemma}\label{lem_univ-split-ex-seq}
{\em(i)}\ \
In the situation of \eqref{sec_dir-summand-conj}, the
following conditions are equivalent :
\begin{enumerate}
\alphaenu
\item
The sequence $\Sigma$ is universally exact.
\item
The sequence $X\otimes_A\Sigma$ is short exact for every
finitely presented $A$-module $X$.
\item
The sequence $\Sigma^\vee:=\Hom_\Z(\Sigma,\Q/\Z)$ splits.
\item
$\Sigma^\vee$ is universally exact.
\item
The sequence $\Hom_A(X,\Sigma)$ is exact for every finitely
presented $A$-module $X$.
\item
$\Sigma$ is the colimit of a filtered system of split short
exact sequences of $A$-modules.
\end{enumerate}

{\em(ii)}\ \
If moreover $M''$ is finitely presented, then
conditions {\em(a)--(f)} are also equivalent to :
\begin{itemize}
\item[(g)]
$\Sigma$ is a split short exact sequence of $A$-modules.
\end{itemize}

{\em(iii)}\ \
If furthermore, $A$ is local and noetherian, and $M,M',M''$
are finitely generated, then conditions {\em(a)--(g)} are
also equivalent to :
\begin{itemize}
\item[(h)]
The sequence $X\otimes_A\Sigma$ is short exact for every
$A$-module $X$ of finite length.
\end{itemize}
\end{lemma}
\begin{proof} Obviusly (f)$\Rightarrow$(a)$\Rightarrow$(b).
To see that (b)$\Rightarrow$(a), write any given $A$-module
$X$ as the filtered colimit of a system
$(X_\lambda~|~\lambda\in\Lambda)$ of finitely presented
$A$-modules; then $H^i(X\otimes_A\Sigma)$ is the colimit
of the induced filtered system of $A$-modules
$(H^i(X_\lambda\otimes_A\Sigma)~|~\lambda\in\Lambda)$, for every
$i\in\Z$. But (b) says that each $X_\lambda\otimes_A\Sigma$ is
exact, whence (a).

(a)$\Rightarrow$(c): Clearly $\Sigma^\vee$ is short exact, since
$\Q/\Z$ is an injective $\Z$-module. For every $A$-module $X$,
recall that $X^\vee:=\Hom_\Z(X,\Q/\Z)$ is naturally an $A$-module
(see \eqref{subsec_Pontryagin}); notice that for every pair of
$A$-modules $X,Y$ we have a natural isomorphism of $A$-modules :
$$
\omega_{X,Y}:\Hom_A(X^\vee,Y^\vee)\isom(X^\vee\otimes_AY)^\vee
\qquad
\phi\mapsto(h\otimes y\mapsto\phi(h)(y)).
$$
There follows a commutative diagram :
$$
\xymatrix@C+60pt{ \Hom_A(M'^\vee,M^\vee)
\ar[r]^-{\Hom_A(M'^\vee,f^\vee)} \ar[d]_{\omega_{M',M}} &
\Hom_A(M'^\vee,M'^\vee) \ar[d]^{\omega_{M',M'}} \\
(M'^\vee\otimes_AM)^\vee \ar[r]^-{(M'^\vee\otimes_Af)^\vee} &
(M'^\vee\otimes_AM')^\vee.
}$$
By assumption, $M'^\vee\otimes_Af$ is injective, hence
$(M'^\vee\otimes_Af)^\vee$ is surjective, so the same holds
for $\Hom_A(M'^\vee,f^\vee)$, and (c) follows easily.

Obviously (c)$\Rightarrow$(d). In order to show that
(d)$\Rightarrow$(e), let $X$ be any finitely presented
$A$-module; by assumption, $X\otimes_Ag^\vee$ is an injection.
On the other hand, for every pair of $A$-modules $X,Y$ we
have a natural $A$-linear map
$$
\sigma_{X,Y}:X\otimes_AY^\vee\to\Hom_A(X,Y)^\vee
\qquad
x\otimes t\mapsto(h\mapsto t\circ h(x))
$$
and it is easily seen that $\sigma_{X,Y}$ is an isomorphism
if $X$ is finitely presented (details left to the reader).
By naturality of $\sigma_{X,\bullet}$, it follows that
$\Hom_A(X,g^\vee)^\vee$ is injective, hence $\Hom_A(X,g^\vee)$
is surjective, whence (e).

(e)$\Rightarrow$(f): Write $M''$ as the colimit of a filtered
system $(M''_\lambda~|~\lambda\in\Lambda)$ of finitely presented
$A$-modules, and let $(j_\lambda:M_\lambda\to M~|~\lambda\in\Lambda)$
be the universal co-cone; for every $\lambda$, set
$\Sigma_\lambda:=j_\lambda^*\Sigma$. Then $\Sigma$ is the
colimit of the resulting system of short exact sequences
$(\Sigma_\lambda~|~\lambda\in\Lambda)$, and we are reduced
to checking that each $\Sigma_\lambda$ is split. However,
(e) implies that for every $\lambda$ there exists an $A$-linear
map $h_\lambda:M''_\lambda\to M$ such that $g\circ h_\lambda=j_\lambda$;
the assertion follows straightforwardly.

Next, if $M''$ is finitely presented, (e) implies that
$\Hom_A(M'',g)$ is surjective, whence (g).

(iii): Let $(A,\fm)$ be a local and noetherian ring, and
let $X$ be any finitely presented $A$-module; if (h) holds,
then for every $n\in\N$ the map $f_n:=X/\fm^nX\otimes_Af$ is
injective. Let $t\in\Ker\,(X\otimes_Af)\setminus\{0\}$;
since the $\fm$-adic topology of $X\otimes_AM'$ is separated
(\cite[Th.8.10(i)]{Mat}), there exists $n\in\N$ such that
the image $\bar t$ of $t$ in $X/\fm^nX\otimes_AM'$ does not
vanish, so $f_n(\bar t)\neq 0$ by assumption, and therefore
$f(t)\neq 0$, a contradiction. So $X\otimes_Af$ is injective
for every such $X$, whence (b).
\end{proof}

\sset\subsubsection{}\label{subsec_enter-Sigma}
Let $\cM$ be a regular log-structure on either the Zariski
or the \'etale site of a noetherian and affine scheme $X$;
say that $X=\Spec\,A$ for a noetherian ring $A$. Let also
$f:A\to B$ be a finite injective ring homomorphism; we want
to show that there exists an $A$-linear map $g:B\to A$ such
that $g\circ f=\one_A$, or equivalently, that the short
exact sequence of $A$-modules
$$
\Sigma_f \qquad : \qquad
0\to A\xrightarrow{f}B\to Q:=\Coker\,f\to 0
\qquad\qquad\qquad
$$
splits. The latter holds if and only if the class
$[\Sigma_f]\in\Ext_A^1(Q,A)$ of $\Sigma_f$ vanishes. Let us recall:

\begin{lemma}\label{lem_Ext-flat-base-change}
{\em (i)}\ \
Let $A$ be a noetherian ring, and $M,N$ two $A$-modules, with
$M$ of finite type; then every flat $A$-algebra $A'$ induces
an isomorphism :
$$
A'\otimes_A\Ext^1_A(M,N)\isom\Ext^1_{A'}(A'\otimes_AM,A'\otimes_AN).
$$

{\em (ii)}\ \
In the situation of \eqref{sec_dir-summand-conj}, suppose
that $A$ is noetherian and $M''$ is of finite type. Then the
following conditions are equivalent :
\begin{enumerate}
\alphaenu
\item
$\Sigma$ is a split short exact sequence of $A$-modules.
\item
For every $\fp\in\Spec\,A$ there exists a faithfully flat
$A_\fp$-algebra $A'$ such that $A'\otimes_A\Sigma$ is a split
short exact sequence of $A'$-module.
\end{enumerate}
\end{lemma}
\begin{proof}(i): Pick a resolution $P^\bullet\to M$ by projective
$A$-modules of finite type; then $\Ext^1_A(M,N)$ is computed
by $H^1\Hom_A(P^\bullet,N)$, and likewise
$\Ext^1_{A'}(A'\otimes_AM,A'\otimes_AN)$ is computed by
$H^1\Hom_{A'}(A'\otimes_AP^\bullet,A'\otimes_AN)$. But since
$A$ is noetherian and $P^i$ is an $A$-module of finite type
for every $i$, the induced morphism of complexes
$$
A'\otimes_A\Hom_A(P^\bullet,N)\to
\Hom_{A'}(A'\otimes_AP^\bullet,A'\otimes_AN)
$$
is an isomorphism, whence the assertion.

(ii): It is easily seen that the class
$[A'\otimes_A\Sigma]\in\Ext^1_{A'}(A'\otimes_AQ,A')$ corresponds
-- under the isomorphism of (i) -- to
$1\otimes[\Sigma]\in A'\otimes_A\Ext^1_A(Q,A')$. The assertion
follows immediately.
\end{proof}

In view of lemma \ref{lem_Ext-flat-base-change}(ii) we are
reduced to checking that $f_\fp:A_\fp\to B_\fp$ admits an
$A_\fp$-linear section for every $\fp\in\Spec\,A$; hence {\em we
may assume that $(A,\fm)$ is a local noetherian ring}. Likewise,
let $A^\sh$ be a strict henselization of $A$ (at the closed point
$\fm$); since the structure map $A\to A^\sh$ is faithfully flat,
lemma \ref{lem_Ext-flat-base-change}(ii) implies that $[\Sigma_f]=0$
if and only if $[A^\sh\otimes_A\Sigma_f]=0$ in
$\Ext^1_{A^\sh}(A^\sh\otimes_AQ,A^\sh)$; also, $A^\sh$ is noetherian,
by \cite[Ch.IV, Prop.18.8.8(iv)]{EGA4}. Thus, it suffices to check
that $A^\sh\otimes_Af$ admits an $A^\sh$-linear section, and
therefore {\em we may assume that $\cM$ admits a chart
$\beta:P\to A$ that is sharp at the point $\fm$}.

\sset\subsubsection{}\label{subsec_jeremy}
Arguing likewise, we may replace $A$ by its $\fm$-adic completion,
and suppose therefore that $A$ is a complete local noetherian
ring. Moreover, $A$ is normal and Cohen-Macaulay, by corollary
\ref{cor_normal-and-CM}. Furthermore, the proof of theorem
\ref{th_charact-log-regular} show that there exists $r\in\N$
and a surjective ring homomorphism
$$
\phi:V[[Q]]\to A
\qquad
\text{where $Q:=P\times\N^{\oplus r}$}
$$
and where :
\begin{itemize}
\item
$(V,\fm_V)$ is a {\em coefficient ring} of $A$, {\em i.e.} either
a field, or a complete discrete valuation ring with $A=\phi(V)+\fm$
and such that $\fm_V=pV$ for a prime integer $p$
\item
The kernel of $\phi$ is trivial if $V$ is a field, and otherwise
it is principal, generated by a power series
$\sum_{q\in Q}a_q\cdot q$ with $a_0\in\fm_V\setminus\fm_V^2$
\item
The chart $\beta$ is the composition of $\phi$ with the
obvious inclusion map $P\to Q\to V[[Q]]$.
\end{itemize}
Let $k'$ be an algebraically closed field extension of the residue
field $k$ of $A$; by \cite[Th.29.1]{Mat} there exists a flat ring
homomorphism $\psi:V\to V'$, where $V'=k'$ if $V=k$, and otherwise
$V'$ is a complete discrete valuation ring with maximal ideal $pV'$
and residue field $k'$. The induced ring extension
$\psi[Q]:V[Q]\to V'[Q]$ is again faithfully flat, hence it induces
a faithfully flat ring homomorphism $V[[Q]]\to V'[[Q]]$
(\cite[Th.22.4]{Mat}). Set $A':=V'[[Q]]\otimes_{V[[Q]]}A$; the
induced map $A\to A'$ is faithfully flat, and $A'$ is the quotient
of $V'[[Q]]$ by the principal ideal generated by the image of
$\Ker\,\phi$. Hence, $A'$ is a complete local noetherian ring, and
by theorems \ref{th_charact-log-regular} and \ref{th_reg-generizes},
the induced map of monoids $Q\to A'$ is a chart for a regular log
structure on $\Spec\,A'$. By lemma \ref{lem_Ext-flat-base-change}(ii),
we may then replace $A,B,P$ by respectively $A',A'\otimes_AB,Q$, and
{\em assume that $k$ is algebraically closed and $\fm=\fm_PA$} (recall
that $\fm_P=P\setminus\{0\}$ : see \eqref{subsec_sepc-of-monoid}).

\sset\subsubsection{}\label{subsec_reduce-to-domain}
Let $\fp_1,\dots,\fp_k$ be the finitely many minimal prime ideals
of $B$; since $A$ is reduced, $f(A)\cap\fp_1\cap\cdots\cap\fp_k=0$,
hence the composition of $f$ with the natural map
$B\to B/\fp_1\times\cdots\times B/\fp_k$ is still finite and injective,
and its kernel contains the ideal $\prod_{i=1}^kf^{-1}\fp_i$. Then
there exists $i\in\{1,\dots,k\}$ such that $f^{-1}\fp_i=0$, so that
the composition $\bar f_i:A\to B/\fp_i$ of $f$ with the projection
$B\to B/\fp_i$ is still injective; clearly it suffices to check
that $\bar f_i$ admits an $A$-linear section, so {\em we may assume
that $B$ is a domain}.

Let $K$ be the field of fractions of $A$, set $B_K:=K\otimes_AB$,
let $d:=\dim_KB_K$, and denote by $\Tr_{B_K/K}:B_K\to K$ the trace
form; since $A$ is integrally closed in $K$, the composition of
$\Tr_{B_K/K}$ with the localization $B\to B_K$ factors through an
$A$-linear map (lemma \ref{lem_axiomatize-familiar}(i))
$$
\Tr_{B/A}:B\to A.
$$
Then $\Tr_{B/A}\circ f=d\cdot\one_A$. In particular, if $A$
is a $\Q$-algebra, the map $d^{-1}\cdot\Tr_{B/A}$ is a section
of $f$ as sought. Thus, {\em we may assume that $k$ has
characteristic $p>0$}. Let us set
$$
P^{(n)}:=\{\gamma\in P_\Q~|~\gamma^{p^n}\in P\}
\qquad
P^{(\infty)}:=\textstyle{\bigcup_{n\in\N}P^{(n)}}
\qquad\text{and}\qquad
\Gamma:=(P^{(\infty)}/P)^\gp.
$$
Then $P^{(\infty)}$ is naturally a $\Gamma$-graded monoid.
More precisely, for every $\gamma\in(P^{(\infty)})^\gp$, denote
by $[\gamma]\in\Gamma$ the class of $\gamma$; then we have :
$$
\gr_{[\gamma]}P^{(\infty)}=P^{(\infty)}\cap\gamma P^\gp
\qquad
\text{for every $\gamma\in(P^{(\infty)})^\gp$}.
$$

\begin{theorem}\label{th_dir-summand-conj-F_p-alg}
In the situation of \eqref{subsec_enter-Sigma}, suppose that
$A$ is an $\F_p$-algebra. Then $f$ admits an $A$-linear section.
\end{theorem}
\begin{proof} After the foregoing reductions, we may assume
that $A=k[[P]]$ for a sharp, fine and saturated monoid $P$,
and that $B$ is a domain. For every $n\in\N$, set
$$
A_n:=k[[P^{(n)}]]
\qquad
B_n:=A_n\otimes_AB
\qquad\text{and}\qquad
A_\infty:=\textstyle{\bigcup_{n\in\N}A_n}.
$$
Notice that $A_0$ is a direct summand of the $A_0$-module
$A_n$, hence the inclusion map $i_n:A\to A_n$ has a section
$\sigma_n:A_n\to A$ for every $n\in\N$; we are then reduced
to checking that the induced map $f_n:A_n\to B_n$ has an
$A_n$-linear section $s_n:B_n\to A_n$, for some $n\in\N$ :
indeed, in this case let $j_n:B\to B_n$ be the induced
inclusion map, and set $s:=\sigma_n\circ s_n\circ j_n:B\to A$.
We get:
$$
s\circ f=\sigma_n\circ s_n\circ j_n\circ f=
\sigma_n\circ s_n\circ f_n\circ i_n=\sigma_n\circ i_n=\one_A
$$
as required. On the other hand, let $K_n$ be the field
of fractions of $A_n$, and denote by $B'_n$ the maximal
reduced quotient of $B_n$, for every $n\in\N$; since $A_n$
is reduced, the induced map $f'_n:A_n\to B'_n$ is still
injective, so it suffices to exhibit a section for $f'_n$.
Moreover, $E_n:=K_n\otimes_{A_n}B'_n$ is a finite field extension
of $K_n$ (since $K_n$ is a purely inseparable field extension
of $K$), and the inclusion map $K_n\to K_{n+1}$ induces a
field extension $E_n\to E_{n+1}$ for every such $n$. Notice
also that $K_\infty:=\bigcup_{n\in\N}K_n$ is a perfect field; hence
$E_\infty:=\bigcup_{n\in\N}E_n$ is the maximal reduced quotient
of $K_\infty\otimes_AB$, and is a finite {\em separable} field
extensions of $K_\infty$. It follows easily that for $n\in\N$
sufficiently large, $E_n$ is already a {\em separable and
finite} field extension of $K_n$.
Summing up, we may replace $A$ and $B$ by $A_n$ and $B'_n$
for $n$ large enough, and {\em assume that $f$ is generically
\'etale, i.e.} there exists $g\in A\setminus\{0\}$ such that
$A[1/g]\otimes_Af$ is a finite, injective and \'etale map.
Then, notice that $A_\infty$ is a perfect $\F_p$-algebra, and
let $\fn\subset A_\infty$ be the radical of the ideal $A_\infty g$;
let also $C_\infty$ be the integral closure of $A_\infty$ in
$A_\infty[1/g]\otimes_AB$. By theorems \ref{th_perfect-purity}
and \ref{th-alm-purity-form-perfectoid}(ii), and proposition
\ref{prop_almost-pure-crit}, the $(A_\infty,\fn)^a$-algebra
$C_\infty^a$ is \'etale, faithfully flat and of almost finite
rank. Now, let $X$ be any $A$-module of finite length; by
virtue of lemma \ref{lem_univ-split-ex-seq}(ii,iii), it suffices
to check that the map $X\otimes_Af:X\to X\otimes_AB$ is
injective. However, since $A$ is a direct summand of the
$A$-module $A_\infty$, the induced map $h:X\to X\otimes_AA_\infty$
is injective, and since $C^a_\infty$ is a faithfully flat
$(A_\infty,\fn)^a$-algebra, the kernel of the induced map
$l:X\otimes_AA_\infty\to X\otimes_AC_\infty$ is annihilated by
$\fn$. Clearly $l\circ h$ factors through $X\otimes_Af$; hence :
$$
\fn\cdot\Ker(X\otimes_Af)=0
$$
(we regard $\Ker(X\otimes_Af)$ as a subset of $\Ker\,l$, via
$h$). Now, as $P^{(\infty)}$ is $\Gamma$-graded, $A_\infty$ is
naturally a $\Gamma$-graded $A$-algebra, and $X\otimes_AA_\infty$
is a $\Gamma$-graded $(A_\infty,\gr_\bullet A_\infty)$-module (see
definition \ref{def_Gamma-graded-algs}(ii)); by definition :
$$
\gr_{[\gamma]}A_\infty:=A\otimes_{\Z[P]}\Z[\gr_{[\gamma]}P^{(\infty)}]
\qquad
\text{for every $\gamma\in P^{(\infty)}$}.
$$
According to proposition \ref{prop_reflex-dual}(ii), there
exist $\lambda_1,\dots,\lambda_k\in P^\vee:=\Hom_\Mnd(P,\N)$
such that $\R\lambda_1,\dots,\R\lambda_k$ are the extremal
rays of $P^\wedge_\R$ (see \eqref{subsec_extremal}).

\begin{claim}\label{cl_good-gammas}
If $\gamma\in P^{(\infty)}$, and $0\leq\lambda_i(\gamma)<1$ for
every $i=1,\dots,k$, then $\gr_{[\gamma]}P^{(\infty)}=\gamma P$.
\end{claim}
\begin{pfclaim} By proposition \ref{prop_reflex-dual}(iv),
we have
$$
P^{(\infty)}=\{\gamma\in(P^{(\infty)})^\gp~|~
\lambda_i(\gamma)\geq 0\text{ for every }i=1,\dots,k\}.
$$
Hence $\gr_{[\gamma]}P^{(\infty)}$ is the subset of all
$\gamma\cdot\beta$, with $\beta\in P^\gp$ and
$\lambda_i(\gamma\cdot\beta)\geq 0$ for $i=1,\dots,k$. But
$\lambda_i(\beta)\in\Z$ for every such $i$ and $\beta$; our
assumption on $\lambda_i(\gamma)$ then implies that
$\lambda_i(\beta)\in\N$ for every such $i$ and $\beta$, whence
$\beta\in P$ (again, by proposition \ref{prop_reflex-dual}(iv)),
and the claim follows.
\end{pfclaim}

We have $g=\sum_{\gamma\in P}g_\gamma\cdot\gamma$ for a system
$(g_\gamma~|~\gamma\in P)$ of elements of $k$. Pick $\beta\in P$
with $g_\beta\neq 0$, and $N\in\N$ large enough, so that
$\lambda_i(\beta^{1/p^N})=p^{-N}\lambda_i(\beta)<1$ for
$i=1,\dots,k$; by claim \ref{cl_good-gammas}, we then have
$\gr_{[\beta^{1/p^N}]}P^{(\infty)}=\beta^{1/p^N}P$. Clearly
$g^{1/p^N}=\sum_{\gamma\in P}g^{1/p^N}_\gamma\cdot\gamma^{1/p^N}$; let
us also write $g^{1/p^N}=\sum_{[\gamma]\in\Gamma}g_{[\gamma]}$, with
$g_{[\gamma]}\in\gr_{[\gamma]}A_\infty$ for every $[\gamma]\in\Gamma$.
It follows that
$g_{[\beta^{1/p^N}]}=g^{1/p^N}_\beta\cdot\beta^{1/p^N}\cdot(1+a)$ for
some $a\in\fm$. Let now $x\in\Ker\,X\otimes_Af$; we have
$x\in\gr_0(X\otimes_AA_\infty)$, and $g^{1/p^N}x=0$; then
$g_{[\gamma]}x=0$ for every $[\gamma]\in\Gamma$, and especially :
$$
x\otimes g^{1/p^N}_\beta\cdot\beta^{1/p^N}=0
\quad\text{in}\quad
\gr_{[\beta^{1/p^N}]}(X\otimes_AA_\infty)=
X\otimes_AA_{[\beta^{1/p^N}]}=X\otimes_{\Z[P]}\Z[\beta^{1/p^N}P]
$$
and since $g_\beta\neq 0$, we have $x\otimes \beta^{1/p^N}=0$ in
$X\otimes_{\Z[P]}\Z[\beta^{1/p^N}P]$, so finally $x=0$, as required.
\end{proof}

\sset\subsubsection{}\label{subsec_mix-char-case}
By theorem \ref{th_dir-summand-conj-F_p-alg} and the foregoing
reductions, we may assume that
$$
A=V[[P]]/(\theta)
$$ 
for a fine, sharp and saturated monoid $P\neq 0$, a complete
discrete valuation ring $(V,\fm_V)$ with $\fm_V=pV$, whose
residue field $k:=V/pV$ is algebraically closed, and a power
series $\theta:=\sum_{\gamma\in P}a_\gamma\cdot\gamma$ with
$a_0\in pV\setminus p^2V$. Then $\fm_PA=\fm$, the maximal
ideal of $A$. We set
$$
A_n:=A\otimes_{\Z[P]}\Z[P^{(n)}]
\qquad\text{and}\qquad
B_n:=A_n\otimes_AB
\qquad
\text{for every $n\in\N$}.
$$
Notice that the induced map $A_n\to A_{n+1}$ is injective for
every $n\in\N$, since $\Z[P^{(n)}]$ is a direct summand of
the $\Z[P^{(n)}]$-module $\Z[P^{(n+1)}]$; hence, let as well
$A_\infty:=\bigcup_{n\in\N}A_n$. Moreover, notice that
$A/\fm_PA=V/a_0V=k$, and $\Omega^1_{k/\Z}=0$, since $k$ is
perfect; then, the integer $r$ of \eqref{subsec_more-notation}
vanishes in the current situation, hence theorem
\ref{th_not-an-F_p}(iv) tells us that $A_\infty$ {\em is a
formal perfectoid ring, for its $p$-adic topology}, and clearly
it is also naturally $\Gamma$-graded.

Let $\lambda_1,\dots,\lambda_k\in P^\vee:=\Hom_\Mnd(P,\N)$
such that $\R\lambda_1,\dots,\R\lambda_k$ are the extremal
rays of $P^\vee_\R$ (see the proof of theorem
\ref{th_dir-summand-conj-F_p-alg}), and
$\beta_1,\dots,\beta_l\in P$ a finite system of generators
of the $\Q_+$-module $\Q_+\otimes_\N P$ (see
\eqref{subsec_adjunct-scoppia}); set
$$
\lambda:=\lambda_1+\cdots+\lambda_k
\qquad
\beta:=\beta_1\cdots\beta_l
\qquad\text{and}\qquad
L:=\max(\lambda(\beta_1),\dots,\lambda(\beta_l)).
$$

\begin{lemma}\label{lem_explicit-epsilon}
In the situation of \eqref{subsec_mix-char-case}, let
$J\subset A_\infty$ be a $\Gamma$-graded ideal,
$g\in A_\infty$, and let $e,n\in\N$ such that
$$
g\in\fm^eA_\infty\setminus\fm^{e+1}A_\infty
\qquad
g\in J^{p^n}
\qquad\text{and}\qquad
1\geq\eps(e,n):=(eL+\lambda(\beta))/p^n.
$$
Then there exists $\gamma\in P^{(\infty)}$ such that :
$$
\gr_{[\gamma]}J=\gr_{[\gamma]}A_\infty
\qquad\text{and}\qquad
\lambda(\gamma)<\eps(e,n).
$$
\end{lemma}
\begin{proof} We may write $g=g_1+\cdots+g_r$ for a finite
sequence $g_1,\dots,g_r\in J^{p^n}$ such that each $g_i$ is
in turn of the form $g_i=\prod_{j=1}^{p^n}a_{ij}\otimes\gamma_{ij}$
with :
$$
a_{i1},\dots,a_{i,p^n}\in A
\qquad
\gamma_{i1},\dots,\gamma_{i,p^n}\in P^{(\infty)}
\qquad
a_{ij}\otimes\gamma_{ij}\in\gr_{[\gamma_{ij}]}J
\quad\text{for every $j=1,\dots,p^n$}.
$$
It is clear that there exists $e'\leq e$ and $i\in\{1,\dots,r\}$
such that $g_i\in\fm^{e'}A_\infty\setminus\fm^{e'+1}A_\infty$. Up to
reordering the $p^n$ factors of $g_i$, we may then assume that
there exists $c\in\{0,\dots,p^n\}$ such that
$\lambda(\gamma_{ij})<\eps(e,n)$ for $j=1,\dots,c$, and
$\lambda(\gamma_{ij})\geq\eps(e,n)$ for $j=c+1,\dots,p^n$.
Suppose now that the lemma fails; then we must have $a_{ij}\in\fm$
for every $i=1,\dots,c$, and therefore
\set\begin{equation}\label{eq_hopefully-ok}
\prod_{j=c+1}^{p^n}a_{ij}\otimes\gamma_{ij}\notin\fm^{e'-c+1}A_\infty.
\end{equation}
For every $j=c+1,\dots,p^n$, write
$\gamma_{ij}=\prod_{s=1}^l\beta_s^{d_{ijs}}$ for certain
$d_{ij1},\dots,d_{ijl}\in\Q_+$; then there exist unique
$d'_{is}\in\N$ and $d''_{is}\in[0,1[$ with
$$
\sum_{j=c+1}^{p^n}d_{ijs}=d'_{is}+d''_{is}
\qquad
\text{for every $s=1,\dots,l$}
$$
and \eqref{eq_hopefully-ok} then implies that
$\sum_{s=1}^ld'_{is}\leq e'-c$, whence :
$$
\sum_{j=c+1}^{p^n}\lambda(\gamma_{ij})<
(e'-c)\cdot L+\lambda(\beta')\leq p^n\cdot\eps(e,n)-cL.
$$
However,
$\sum_{j=c+1}^{p^n}\lambda(\gamma_{ij})\geq(p^n-c)\cdot\eps(e,n)$,
and since $\eps(e,n)\leq 1\leq L$, this is absurd.
\end{proof}

\begin{theorem}\label{th_log-direct-summand-conj}
In the situation of \eqref{subsec_enter-Sigma}, the map $f$
admits an $A$-linear section.
\end{theorem}
\begin{proof} In view of the foregoing reductions, we may assume
that $A$ and $A_\infty$ are as in \eqref{subsec_mix-char-case} and
$B$ is a domain; recall also that $B[1/g]$ is an \'etale and
faithfully flat $A[1/g]$-algebra. Let :
$$
A'_n:=A_\infty[X^{1/p^n}]/(X-g)
\qquad\text{for every $n\in\N$}
\qquad\text{and}\qquad
A'_\infty:=\textstyle{\bigcup_{n\in\N}A'_n}.
$$
Let also $D$ be the $p$-integral closure of $A'_\infty$ in
$A'_\infty[1/p]$, and endow $D$ with its $p$-adic topology;
by theorem \ref{th_Zimbabwe}, we know that $D$ is a faithfully
flat $A_\infty$-algebra, and that $D$ is a formal perfectoid ring.
Denote by $D^\wedge$ the $p$-adic completion of $D$, and set
$B':=D^\wedge\otimes_AB[1/g]$. Hence, $B'$ is a finite \'etale
$D^\wedge[1/g]$-algebra. Let now $\fn\subset D^\wedge$ be the
ideal generated by $(g^{1/p^n}~|~n\in\N)$; then $(D^\wedge,\fn)$
is a basic setup, and we set
$D^\wedge_1:=(D^\wedge,\fn)^a_*\subset D^\wedge[1/g]$ (see lemma
\ref{lem_identify-A^a_0*}); lastly, let $B'_1$ be the integral
closure of $D^\wedge_1$ in $B'$. By theorem
\ref{th_perf-Abhyankar}(ii,iii), for every $n\in\N$, the
$(D/p^nD)^a$-algebra $(B'_1/p^nB'_1)^a$ is faithfully flat and
\'etale of finite rank, for the almost structure given by the
basic setup $(D^\wedge_1,\fn)$. After these preliminaries, let
$X$ be any $A$-module of finite length; by lemma
\ref{lem_univ-split-ex-seq}(ii,iii), it suffices to show that
the map $X\otimes_Af$ is injective. However, let
$x\in\Ker(X\otimes_Af)$, and pick $m\in\N$ large enough, so that
$p^mX=0$. Since $A_\infty$ is a $\Gamma$-graded $A$-algebra with
$\gr_0A_\infty=A$, the induced map
$\phi:X\to X_\infty:=X\otimes_AA_\infty/p^mA_\infty$ is injective,
and $X_\infty$ is naturally a $\Gamma$-graded
$(A_\infty,\gr_\bullet)$-module, with $\gr_0X_\infty=X$. Especially,
the annihilator $J$ of $\phi(x)=x\otimes 1\in\gr_0X_\infty$ is a
$\Gamma$-graded ideal of $A_\infty$. Next, since $D$ is a
faithfully flat $A_\infty$-algebra, the induced map
$\phi':X_\infty\to X_\infty\otimes_{A_\infty}D/p^mD$ is again
injective, and $\Ann_D(\phi'(x\otimes 1))=JD$. Lastly, since
$(B'_1/p^mB'_1)^a$ is a faithfully flat $(D/p^mD)^a$-algebra,
the kernel of the induced map
$\phi'':X_\infty\otimes_{A_\infty}D/p^mD\to
X_\infty\otimes_{A_\infty}B'_1/p^mB'_1$ is annihilated by $\fn$.
Clearly $\phi''\circ\phi'\circ\phi$ factors through
$X\otimes_Af$, hence $\phi'(x\otimes 1)\in\Ker\,\phi''$,
and consequently :
$$
\fn\subset JD.
$$
Thus, $g^{1/p^n}\in JD\cap A_\infty=J$ for every $n\in\N$
(\cite[Th.7.5(ii)]{Mat}), {\em i.e.} $g\in J^{p^n}$ for
every $n\in\N$. Notice that the $\fm$-adic topology is
separated on $A_\infty$; hence, let $e\in\N$ such that
$g\in\fm^eA_\infty\setminus\fm^{e+1}A_\infty$; by lemma
\ref{lem_explicit-epsilon}, for every $n\in\N$ there
exists $\gamma_n\in P^{(\infty)}$ such that
$\lambda(\gamma_n)<\eps(e,n)$ and
$\gr_{[\gamma_n]}J=\gr_{[\gamma_n]}A_\infty$. By claim
\ref{cl_good-gammas}, it follows that for every sufficiently
large $n\in\N$, we have
$\gr_{[\gamma_n]}J=A\otimes_{\Z[P]}\Z[\gamma_n+P]$, a free
$A$-module of rank one, with generator $1\otimes\gamma_n$.
Finally, we get $x\otimes\gamma_n=0$ in
$\gr_{[\gamma_n]}X_\infty=X\otimes_{\Z[P]}\Z[\gamma_n+P]$ for
every such $n$; thus, $x=0$, and the proof is concluded.
\end{proof}

\subsection{Diagonal idempotents of generically \'etale maps}
\label{sec_diagonal-idempotents}
Let $f:A\to B$ be a finite and injective ring homomorphism of
noetherian rings; set $X:=\Spec\,A$, denote by $\cB$ the
quasi-coherent $\cO_X$-algebra associated with $B$, and by
$\mu:\cB\otimes_{\cO_X}\cB\to\cB$ the multiplication morphism. For
every affine open subset $U\subset X$ such that $\cB_{|U}$ is an
\'etale $\cO_U$-algebra, we have a unique {\em diagonal idempotent}
$e_{f,U}\in\cB(U)\otimes_{\cO_X(U)}\cB(U)$ such that :
$$
\mu_U(e_{f,U})=1
\quad\text{and}\quad
x\cdot e_{f,U}=0
\quad\text{for every
$x\in\Ker\,(\cB(U)\otimes_{\cO_X(U)}\cB(U)\xrightarrow{\mu_U}\cB(U))$}.
$$
We suppose moreover that $f$ is {\em generically \'etale}, {\em i.e.}
that the union $U_f\subset X$ of all such affine open subsets is
dense in $X$; the uniqueness property of the diagonal idempotents
then implies that there exists a unique section
$$
e_f\in\Gamma(U_f,\cB\otimes_{\cO_X}\cB)
$$
such that $e_{f|U}=e_{f,U}$ for every affine open subset $U\subset U_f$.
Moreover, $\cB_{|U_f}$ is a locally free $\cO_{U_f}$-module of finite
rank, hence we have a well defined $\cO_{U_f}$-linear {\em trace map}
$$
\tr_{\cB/\cO_X}:\cB_{|U_f}\to\cO_{U_f}
$$
that assigns to every affine open subset $U\subset U_f$ and every
$b\in\cB(U)$ the trace of the $\cO_X(U)$-linear endomorphism
$b\cdot\one_{\cB(U)}$. The corresponding {\em trace form}
$$
t_{\cB/\cO_X}:=\tr_{\cB/\cO_X}\circ\mu:
(\cB\otimes_{\cO_X}\cB)_{|U_f}\to\cO_{U_f}
$$
is a perfect pairing (\cite[Th.4.1.14]{Ga-Ra}), hence it induces
an $\cO_{U_f}$-linear isomorphism
$$
\omega_{\cB/\cO_X}:\cB_{|U_f}\isom(\cB^\vee)_{|U_f}
\qquad\text{where}\qquad
\cB^\vee:=\cHom_{\cO_X}(\cB,\cO_X).
$$

\begin{remark}\label{rem_perfectifcation}
Let $f:A\to B$ be as in \eqref{sec_diagonal-idempotents}, $p\in\N$
a prime integer, and suppose that $A$ is an $\F_p$-algebra. Let
also $h\in A$ such that the localization $f_h:A_h\to B_h$ of $f$
is \'etale. It is easily seen that the natural map (notation of
\eqref{subsec_perfectification})
$$
B^\perf\otimes_{A^\perf}B^\perf\to(B\otimes_AB)^\perf
$$
is bijective. By assumption, $h^ne_f$ is in the image of the
localization $B\otimes_AB\to B_h\otimes_{A_h}B_h$, for some $n\in\N$.
The image $e_f^\perf$ of $e_f$ in $(B\otimes_AB)^\perf_h$ is the
diagonal idempotent of the \'etale map
$A^\perf_h\otimes_Af:A^\perf_h\to B^\perf_h=A^\perf_h\otimes_AB$.
Since $(e^\perf_f)^p=e^\perf_f$, we conclude that $h^{n/p^k}e^\perf_f$
lies in the image of the localization
$(B\otimes_AB)^\perf\to(B\otimes_AB)^\perf_h$, for every $k\in\N$.
\end{remark}

\begin{definition}
(i)\ \
A {\em log ring} $(A,P,\beta)$ is the datum of a ring $A$, an
integral saturated monoid $P$ whose associated abelian group
$P^\gp$ is {\em torsion-free}, and a morphism of monoids
$\beta:P\to A$ from $P$ to the multiplicative monoid $(A,\cdot)$,
furnishing a chart for a {\em log structure} $\cP$ on the Zariski
site of $\Spec\,A$ (see definition \ref{def_chart}(i)), which we
call {\em the log structure of} $(A,P,\beta)$.

(ii)\ \
We say that the log ring $(A,P,\beta)$ is {\em regular} if $A$
is noetherian, $P$ is fine and saturated, and the log structure
of $(A,P,\beta)$ is {\em regular}. In this case, $P^\gp$ is a
free abelian group of finite rank, and $A$ is a product of
finitely many normal domains (corollary \ref{cor_normal-and-CM}).
\end{definition}

\sset\subsubsection{}\label{subsec_intro-log-rings}
Let $(A,P,\beta)$ be any log ring, $f:A\to B$ a ring homomorphism
as in \eqref{sec_diagonal-idempotents}, and set $X:=\Spec\,A$. We
consider the monoid $P^{(\infty)}$ associated with $P$ as in
\eqref{subsec_reduce-to-domain}, and the abelian group
$\Gamma:=(P^{(\infty)}/P)^\gp$. We will likewise use the monoid
$P_\Q:=\Q_+\otimes_\N P$, and the abelian group
$\Gamma_\Q:=(P_\Q/P)^\gp$, and for every $\gamma\in P_\Q$, we let
$[\gamma]\in\Gamma_\Q$ be the class of $\gamma$. Then the ring
$$
A_\infty:=A\otimes_{\Z[P]}\Z[P_\Q]
$$
is naturally a $\Gamma_\Q$-graded $A$-algebra, with graded summands :
$$
A_{[\gamma]}:=A\otimes_{\Z[P]}\Z[\gr_{[\gamma]}P_\Q]
\qquad
\text{for every $\gamma\in P_\Q$}.
$$
Likewise, set $B_\infty:=A_\infty\otimes_AB$; both $B_\infty$
and its dual $A_\infty$-module
$B_\infty^\vee:=\Hom_{A_\infty}(B_\infty,A_\infty)$ inherit natural
$\Gamma_\Q$-graded structures whose graded summands are respectively
$$
B_{[\gamma]}:=A_{[\gamma]}\otimes_AB
\qquad\text{and}\qquad
B^\vee_{[\gamma]}:=\Hom_A(B,A_{[\gamma]})
\qquad
\text{for every $\gamma\in P_\Q$}.
$$
Let $\cO_{X_\infty}$ and $\cB_\infty$ be the quasi-coherent
$\cO_X$-algebras associated with $A_\infty$ and $B_\infty$.
Hence, $\cB_{\infty|U_f}$ is an \'etale $\cO_{X_\infty|U_f}$-algebra,
whose trace morphism induces the isomorphism
$$
\omega_{\cB_\infty/\cO_{X_\infty}}:=
\cO_{X_\infty|U_f}\otimes_{\cO_{U_f}}\omega_{\cB/\cO_X}:
\cB_{\infty|U_f}\isom\cB^\vee_{\infty|U_f}.
$$
Also, the corresponding diagonal idempotent
$e_{f,\infty}\in\Gamma(U_f,\cB_\infty\otimes_{\cO_{X_\infty}}\cB_\infty)$
is the image of $e_f$, under the natural morphism
\set\begin{equation}\label{eq-tirata-orecchi}
\cB\otimes_{\cO_X}\cB\to\cO_{X_\infty}\otimes_{\cO_X}\cB\otimes_{\cO_X}\cB
\isom\cB_\infty\otimes_{\cO_{X_\infty}}\cB_\infty.
\end{equation}
More precisely, since $U_f$ is quasi-compact, the $A_\infty$-module
$\Gamma(U_f,\cB_\infty\otimes_{\cO_{X_\infty}}\cB_\infty)$ is also
$\Gamma_\Q$-graded, and on the other hand, $\gr_0P_\Q=P$, since $P$ is
saturated; then $\gr_0B_\infty=B$, and the map \eqref{eq-tirata-orecchi}
identifies $\Gamma(U_f,\cB\otimes_{\cO_X}\cB)$ with
$\gr_0\Gamma(U_f,\cB_\infty\otimes_{\cO_{X_\infty}}\cB_\infty)$.

\sset\subsubsection{}\label{subsec_both-data}
Let $\underline f:=(f,P,\beta)$ be the datum of a log ring
$(A,P,\beta)$ and a ring homomorphism $f:A\to B$ as in
\eqref{sec_diagonal-idempotents}; with the notation of
\eqref{subsec_intro-log-rings}, we attach to $\underline f$
the composition
\set\begin{equation}\label{eq_was-theta_f}
B^\vee_\infty\otimes_{A_\infty}B^\vee_\infty\xrightarrow{\rho}
\Gamma(U_f,\cB^\vee_\infty\otimes_{\cO_{X_\infty}}\cB^\vee_\infty)
\xrightarrow{\alpha}\Gamma(U_f,\cB_\infty\otimes_{\cO_{X_\infty}}\cB_\infty)
\end{equation}
where $\rho$ is the restriction map of the quasi-coherent
$\cO_X$-module $\cB^\vee_\infty\otimes_{\cO_{X_\infty}}\cB^\vee_\infty$,
and $\alpha$ is induced by the isomorphism
$\omega^{-1}_{\cB_\infty/\cO_{X_\infty}}\otimes_{\cO_{X_\infty|U_\infty}}
\omega_{\cB_\infty/\cO_{X_\infty}}^{-1}$. Set as well :
$$
\sR(B,\Delta):=\bigoplus_{[\gamma]\in\Delta}
B^\vee_{[\gamma]}\otimes_AB^\vee_{[1/\gamma]}
\qquad
\text{for every subset $\Delta\subset\Gamma_\Q$}.
$$
Since \eqref{eq_was-theta_f} is a morphism of graded
$A_\infty$-modules, it restricts to an $A$-linear map
$$
\theta_{\underline f}:\sR(B,\Gamma_\Q)\to
\gr_0\Gamma(U_f,\cB_\infty\otimes_{\cO_{X_\infty}}\cB_\infty)\isom
\Gamma(U_f,\cB\otimes_{\cO_X}\cB).
$$

\begin{remark}\label{rem_functoriality-of-theta_f}
(i)\ \
In the situation of \eqref{subsec_both-data}, let $g:A\to A'$
be a ring homomorphism; set $B':=A'\otimes_AB$, let $f':A'\to B'$
be the induced map, and suppose that
$\underline f':=(f',P,g\circ\beta)$ is also a datum as in
\eqref{subsec_both-data}. Set $X':=\Spec\,A'$, and let $h:X'\to X$
be the morphism of schemes induced by $g$; also, let $\cB'$ be the
quasi-coherent $\cO_{X'}$-algebra arising from $B'$. For every
$[\gamma]\in\Gamma_\Q$ we get a natural $A$-linear map
\set\begin{equation}\label{eq_happy}
B^\vee_{[\gamma]}\to\Hom_A(B,A'_{[\gamma]})\isom
B'^\vee_{[\gamma]}:=\Hom_{A'}(B',A'_{[\gamma]})
\end{equation}
and a direct inspection of the definitions yields a
commutative diagram :
$$
\xymatrix@C+10pt{
\sR(B,\Gamma_\Q) \ar[rr]^-{\theta_{\underline f}} \ar[d] & &
\Gamma(U_f,\cB\otimes_{\cO_X}\cB)
\ar[d] \\
\sR(B',\Gamma_\Q) \ar[r]^-{\theta_{\underline f'}} &
\Gamma(U_{f'},\cB'\otimes_{\cO_{X'}}\cB') \ar[r]^-\rho &
\Gamma(h^{-1}U_f,\cB'\otimes_{\cO_{X'}}\cB')
}$$
whose left (resp. right) vertical arrow is induced by the maps
\eqref{eq_happy} (resp. by the natural morphisms of quasi-coherent
$\cO_X$-algebras $\cO_X\to h_*\cO_{X'}$ and $\cB\to h_*\cB'$), and
where $\rho$ is the restriction map. Also, it is clear that the
right vertical arrow maps $e_f$ to $\rho(e_{f'})$.

(ii)\ \
Let $R$ be any ring, $S$ and $S'$ two finite $R$-algebras whose
underlying $R$-modules are projective (of finite rank). Then
all the trace maps $\tr_{S/R},\tr_{S'/R}$ and $\tr_{S\otimes_RS'/R}$
are well defined, and \cite[Lemma 4.1.3]{Ga-Ra} implies that :
$$
\tr_{S\otimes_RS'/R}=\tr_{S/R}\otimes_R\tr_{S'/R}:S\otimes_RS'\to R.
$$
We deduce the following alternative description of $\theta_{\underline f}$.
First, we have a natural $A_\infty$-linear map
\set\begin{equation}\label{eq_nablus}
B^\vee_\infty\otimes_{A_\infty}B^\vee_\infty\to
(B_\infty\otimes_{A_\infty}B_\infty)^\vee:=\Hom_A(B\otimes_AB,A_\infty)
\end{equation}
given by the rule :
$\phi\otimes\phi'\mapsto(b\otimes b'\mapsto\phi(b)\cdot\phi'(b'))$
for every $\phi,\phi'\in B^\vee_\infty$ and every $b,b'\in B$. Next,
the trace map of the $\cO_{X_\infty}$-algebra
$\cB_\infty\otimes_{\cO_X}\cB_\infty$ induces an $\cO_{X_\infty}$-linear map
$$
\omega_{\cB_\infty\otimes_{\cO_{X_\infty}}\cB_\infty/\cO_{X_\infty}}:
\cB_\infty\otimes_{\cO_X}\cB_\infty\to(\cB_\infty\otimes_{\cO_X}\cB_\infty)^\vee
$$
whose restriction to the open subset $U_f$ is an isomorphism. Then
it is easily seen that \eqref{eq_was-theta_f} equals the composition
of \eqref{eq_nablus} with the map
$$
(B_\infty\otimes_{A_\infty}B_\infty)^\vee\xrightarrow{\rho'}
\Gamma(U_f,(\cB_\infty\otimes_{\cO_X}\cB_\infty)^\vee)
\xrightarrow{\alpha'}\Gamma(U_f,\cB_\infty\otimes_{\cO_X}\cB_\infty)
$$
where $\rho'$ is the restriction map of the quasi-coherent
$\cO_X$-module $(\cB_\infty\otimes_{\cO_X}\cB_\infty)^\vee$, and
$\alpha'$ is induced by the isomorphism
$(\omega_{\cB_\infty\otimes_{\cO_{X_\infty}}\cB_\infty/\cO_{X_\infty}})_{|U_f}$.
Moreover, $(B_\infty\otimes_{A_\infty}B_\infty)^\vee$ is naturally
a $\Gamma_\Q$-graded $A_\infty$-module, with a natural identification :
$$
\gr_0(B_\infty\otimes_{A_\infty}B_\infty)^\vee\isom(B\otimes_AB)^\vee:=
\Hom_A(B\otimes_AB,A)
$$
and \eqref{eq_nablus} is a morphism of $\Gamma_\Q$-graded modules.
It follows that $\theta_{\underline f}$ is the composition
$$
\sR(B,\Gamma_\Q)\xrightarrow{\ \xi_{\underline f}\ }(B\otimes_AB)^\vee
\xrightarrow{\ \xi'_{\underline f}\ }\Gamma(U_f,\cB\otimes_{\cO_X}\cB)
$$
where $\xi_{\underline f}$ is induced by the restriction of
\eqref{eq_nablus}, and $\xi'_{\underline f}$ is the restriction
of $\alpha'\circ\rho'$. In turn, $\xi'_{\underline f}$ can be
decribed as the composition of the restriction map
$(B\otimes_AB)^\vee\to\Gamma(U_f,(\cB\otimes_{\cO_X}\cB)^\vee)$
of the quasi-coherent $\cO_X$-module $(\cB\otimes_{\cO_X}\cB)^\vee$
arising from $(B\otimes_AB)^\vee$, and the isomorphism
$\Gamma(U_f,(\cB\otimes_{\cO_X}\cB)^\vee)\isom
\Gamma(U_f,\cB\otimes_{\cO_X}\cB)$ induced by the isomorphism
$(\cB\otimes_{\cO_X}\cB)^\vee_{|U_f}\isom(\cB\otimes_{\cO_X}\cB)_{U_f}$
deduced as usual from the trace map of $\cB\otimes_{\cO_X}\cB$.

(iii)\ \
With the notation of (ii), notice as well that $\xi'_{\underline f}$
is injective if $A$ is reduced. Indeed, if $\phi:B\otimes_AB\to A$
restricts on $U_f$ to the zero section of $(\cB\otimes_{\cO_X}\cB)^\vee$,
then for every $b,b'\in B$ the element $\phi(b\otimes b')\in A$
restricts on $U_f$ to the zero section of $\cO_X$; since $A$ is
reduced and $U_f$ is dense in $X$, this means that
$\phi(b\otimes b')=0$ for every such $b,b'$, whence the contention.
\end{remark}

\begin{lemma}\label{lem_finite-Sigma}
{\em(i)}\ \
In the situation of \eqref{subsec_both-data}, let
$\Delta\subset\Gamma_\Q$ be any subset. Then there exists
a finite subset $\Sigma\subset\Delta$ such that
$\xi_{\underline f}(\sR(B,\Delta))=\xi_{\underline f}(\sR(B,\Sigma))$
(where $\xi_{\underline f}$ is defined as in remark
{\em\ref{rem_functoriality-of-theta_f}(ii)}). Moreover, $\Sigma$
depends only on $P$ and $\Delta$ (and neither on $\beta$ nor $f$).

{\em(ii)}\ \
Suppose that $P=P_1\times P_2$, where
$P_2=\N^{\oplus r}\times\Z^{\oplus s}$, for some $r,s\in\N$. For $i=1,2$,
let $\Delta_i\subset\Gamma_{i,\Q}:=(P_{i,\Q}/P_i)^\gp$ be any subset,
so that $\Delta_1\times\Delta_2\subset\Gamma_\Q=
\Gamma_{1,\Q}\oplus\Gamma_{2,\Q}$. Then
$$
\xi_{\underline f}(\sR(B,\Delta_1\times\Delta_2))=
\xi_{\underline f}(\sR(B,\Delta_1\times\{0\})).
$$
\end{lemma}
\begin{proof}(i): On the one hand, the multiplication law of
$A_\infty$ induces by restriction $A$-linear maps
$$
\mu_{[\gamma]}:A_{[\gamma]}\otimes_AA_{[1/\gamma]}\to\gr_0A_\infty=A
\qquad
(a\otimes\gamma\lambda)\otimes(a'\otimes\lambda'/\gamma)\mapsto
aa'\cdot\beta(\lambda\lambda')
$$
(notice that if $\gamma,\lambda,\lambda'\in P^\gp_\Q$ and we have
both $\gamma\lambda\in\gr_{[\gamma]}P_\Q$ and
$\lambda'/\gamma\in\gr_{[1/\gamma]}P_\Q$, then
$\lambda\lambda'\in\gr_0P_\Q=P$). On the other hand, we have an
obvious isomorphism of $P$-modules :
$$
Q_\gamma:=P^\gp\cap\gamma^{-1}P_\Q\isom\gr_{[\gamma]}P_\Q
\qquad
\lambda\mapsto\gamma\lambda
\qquad\text{for every $\gamma\in P^\gp_\Q$}.
$$
Especially, if $Q_\gamma=Q_\delta$ for some
$\gamma,\delta\in P^\gp_\gamma$, we deduce an $A$-linear isomorphism
$$
\omega_{\gamma,\delta}:A_{[\gamma]}\isom A_{[\delta]}
\qquad
a\otimes\gamma\lambda\mapsto a\otimes\delta\lambda
\qquad
\text{for every $a\in A$ and $\gamma\lambda\in\gr_{[\gamma]}P_\Q$}
$$
which induces an $A$-linear isomorphism
$$
\omega_{\gamma,\delta*}:B^\vee_{[\gamma]}\isom B^\vee_{[\delta]}
\qquad
(\phi:B\to A_{[\gamma]})\mapsto
(\omega_{\gamma,\delta}\circ\phi:B\to A_{[\delta]}).
$$
Furthermore, if we have both $Q_\gamma=Q_\delta$ and
$Q_{1/\gamma}=Q_{1/\delta}$ for some $\gamma,\delta\in P^\gp_\gamma$,
a direct inspection yields a commutative diagram :
\set\begin{equation}\label{eq_triangular-inspection}
{\diagram\spreaddiagramcolumns{+20pt}
A_{[\gamma]}\otimes_AA_{[1/\gamma]}
\ar[rr]^-{\omega_{\gamma,\delta}\otimes_A\omega_{1/\gamma,1/\delta}}
\ar[rd]_{\mu_{[\gamma]}} & &
A_{[\delta]}\otimes_AA_{[1/\delta]} \ar[ld]^{\mu_{[\delta]}} \\
& A.
\enddiagram}
\end{equation}
For such $\gamma$ and $\delta$, there follows a commutative diagram :
$$
\xymatrix{
B^\vee_{[\gamma]}\otimes_AB^\vee_{[1/\gamma]}
\ar[rr]^-{\omega_{\gamma,\delta*}\otimes_A\omega_{1/\gamma,1/\delta*}}
\ar[rd] & &
B^\vee_{[\delta]}\otimes_AB^\vee_{[1/\delta]} \ar[ld] \\
& (B\otimes_AB)^\vee
}$$
whose downward arrows are the restrictions of $\xi_{\underline f}$ : the
details shall be left to the reader. Pick a section of the projection
$P^\gp_\Q\to\Gamma_\Q$ :
$$
\Gamma_\Q\to P^\gp_\Q\subset P^\gp_\R
\qquad
[\gamma]\mapsto[\gamma]^*
$$
whose image is contained in a {\em bounded} subset of the
finite dimensional $\R$-vector space $P^\gp_\R$. Summing up,
we see that if $Q_{[\gamma]^*}=Q_{[\delta]^*}$ and
$Q_{[1/\gamma]^*}=Q_{[1/\delta]^*}$, then
$\xi_{\underline f}(B^\vee_{[\gamma]}\otimes B^\vee_{[1/\gamma]})=
\xi_{\underline f}(B^\vee_{[\delta]}\otimes B^\vee_{[1/\delta]})$.
To conclude, it suffices now to invoke proposition
\ref{prop_linear-constr-part}(i).

(ii): For every $(\gamma_1,\gamma_2)\in P^\gp_\Q$, let
$[\gamma_i]\in\Gamma_{i,\Q}$ be the class of $\gamma_i$ for
$i=1,2$, and $[\gamma_1,\gamma_2]\in\Gamma_\Q$ the class of
$(\gamma_1,\gamma_2)$; also, for each $x\in\Q$, denote by
$x^\dagger$ the smallest element of $(x+\Z)\cap\Q_+$. It is
easily seen that
$$
\gr_{[\gamma_1,\gamma_2]}P_\Q=
(\gr_{[\gamma_1]}P_{1,\Q})\times(\gamma_2^\dagger P_2)
\qquad
\text{where 
$\gamma_2^\dagger:=(\gamma^\dagger_{2,1},\dots,\gamma^\dagger_{2,r})$}
$$
whence a natural identification :
\set\begin{equation}\label{eq_barbara-top}
B^\vee_{[\gamma_1,0]}\otimes_{\Z[P]}\Z[\gamma_2^\dagger P_2]\isom
B^\vee_{[\gamma_1,\gamma_2]}.
\end{equation}
Notice that for every $x\in\Q$ the rational number
$x^\dagger+(-x)^\dagger$ equals $0$ if $x\in\Z$, and otherwise
it equals $1$. Especially :
$$
\gamma_2^\dagger+(-\gamma_2)^\dagger\in P_2
\qquad
\text{for every $\gamma_2\in P^\gp_{2,\Q}$}.
$$
Therefore the multiplication law of $\Z[P_{2,\Q}]$ induces by
restriction an $A$-linear map :
$$
\mu_{\gamma_2}:
\Z[\gamma_2^\dagger P_2]\otimes_{\Z[P]}\Z[(-\gamma_2)^\dagger P_2]
\to\Z[P_2].
$$
We then get a commutative diagram
$$
\xymatrix{ (B^\vee_{[\gamma_1,0]}\otimes_{\Z[P_2]}\Z[\gamma_2^\dagger P_2])
\otimes_A(B^\vee_{[1/\gamma_1,0]}\otimes_{\Z[P_2]}\Z[(-\gamma_2)^\dagger P_2])
\ar[r] \ar[d]_{\xi_{\underline f,[\gamma_1,0]}\otimes_{\Z[P_2]}
\Z[\gamma_2^\dagger P_2]\otimes_{\Z[P_2]}\Z[(-\gamma_2)^\dagger P_2]} &
B^\vee_{[\gamma_1,\gamma_2]}\otimes_AB^\vee_{[1/\gamma_1,-\gamma_2]}
\ar[d]^{\xi_{\underline f,[\gamma_1,\gamma_2]}} \\
(B\otimes_AB)^\vee\otimes_{\Z[P_2]}
\Z[\gamma_2^\dagger P_2]\otimes_{\Z[P_2]}\Z[(-\gamma_2)^\dagger P_2]
\ar[r] & (B\otimes_AB)^\vee
}$$
where $\xi_{\underline f,[\gamma_1,0]}:
B^\vee_{[\gamma_1,0]}\otimes_AB^\vee_{[1/\gamma_1,0]}\to(B\otimes_AB)^\vee$
is the restriction of $\xi_{\underline f}$, and likewise for
$\xi_{\underline f,[\gamma_1,\gamma_2]}$. Also, the top horizontal
arrow is the isomorphism induced by \eqref{eq_barbara-top}, and the
bottom horizontal arrow is $(B\otimes_AB)^\vee\otimes_{\Z[P_2]}\mu_{\gamma_2}$.
Summing up, we conclude that :
$$
\xi_{\underline f}(B^\vee_{[\gamma_1,\gamma_2]}\otimes_AB^\vee_{[1/\gamma_1,-\gamma_2]})
\subset\xi_{\underline f}(B^\vee_{[\gamma_1,0]}\otimes_AB^\vee_{[1/\gamma_1,0]})
\qquad
\text{for every $[\gamma_1,\gamma_2]\in P^\gp_\Q$}
$$
whence the assertion.
\end{proof}

\begin{example}\label{ex_regular-log-ring}
Let $V$ be a noetherian regular ring, $P$ a fine and saturated
monoid such that $P^\gp$ is torsion-free, $\psi:V[P]\to A$ a
smooth ring homomorphism, and $\beta:P\to A$ the composition
of $\psi$ with the natural inclusion map $\alpha:P\to V[P]$.
Set $S:=\Spec\,V$, $X_0:=\Spec\,V[P]$ and $X:=\Spec\,A$; then
$\alpha$ is a chart for a log structure $\cP_0$ on the Zariski
site of $X_0$, and the induced morphism of log schemes
$(X_0,\cP_0)\to(S,\cO^\times_S)$ is smooth, by virtue of
proposition \ref{prop_toric-smooth}. Since $(S,\cO^\times_S)$ is
trivially a regular log scheme, it follows that the same holds
for $(X_0,\cP_0)$, by theorem \ref{th_smooth-preserve-reg}; then,
again by theorem \ref{th_smooth-preserve-reg} and corollary
\ref{cor_undercover}, the log scheme
$(X,\cP):=X\times_{X_0}(X_0,\cP_0)$ is regular, and $\beta$ is
a chart for $\cP$. Summing up, this shows that the datum
$(A,P,\beta)$ is a regular log ring.
\end{example}

\begin{lemma}\label{lem_from-generic-flatness}
Let $A$ be a domain, $B$ an $A$-algebra of finite type, $C$
a $B$-algebra of finite type, $M$ a $B$-module of finite type,
and $N$ a $C$-module of finite type. Set also $K:=\Frac(A)$.

{\em(i)}\ \
There exists $f\in A\setminus\{0\}$ such that the following
holds for every $A$-algebra $A'$. Set $B':=A'\otimes_AB$,
$C':=A'\otimes_AC$, $M':=B'\otimes_BM$ and $N':=C'\otimes_CN$;
then the natural map
$$
A'\otimes_A\Hom_B(M,N)_f\to\Hom_{B'}(M',N')_f
$$
is an isomorphism.

{\em(ii)}\ \
For every $B$-linear map $\phi:M\to N$ there exists
$f\in A\setminus\{0\}$ such that $(\Coker\,\phi)_f$
is a free $A$-module.

{\em(iii)}\ \
In the situation of {\em(ii)}, let $x\in N$. Then
$x\in\Img\,(K\otimes_A\phi)$ if and only if there exists a dense
subset $U\subset\Spec\,A$ such that
$1\otimes x\in\Img(\kappa(\fp)\otimes_A\phi)$ for every $\fp\in U$.
\end{lemma}
\begin{proof}(i): Pick $n\in\N$ and a $B$-linear surjection
$\psi:B^{\oplus n}\to M$; let also $M'\subset B^{\oplus n}$ be a
finitely generated $B$-submodule such that
$K\otimes_AM'=\Ker\,(K\otimes_A\psi)$. Set $M'':=B^{\oplus n}/M'$,
and denote by $\bar\psi:M''\to M$ the $B$-linear surjection
induced by $\psi$. By \cite[Ch.IV, Lemma 8.9.4.1]{EGAIV-3},
there exists $f\in A\setminus\{0\}$ such that the $A_f$-module
$M''_f$ is flat, and by construction $\Ker(\psi_f:M''_f\to M_f)$
is a torsion $A$-module; then $\Ker\,\psi_f=0$, hence $M_f$ is a
finitely presented $B_f$-module. We may thus replace $A,B,C,M,N$
by $A_f,B_f,C_f,M_f,N_f$, and assume that $M$ is a $B$-module of
finite presentation. Now, let
$$
B^{\oplus m}\xrightarrow{\phi}B^{\oplus n}\to M
$$
be a finite presentation of $M$; there follow short exact
sequences of $B$-modules :
$$
\begin{aligned}
\Sigma \quad:\quad &
0\to\Hom_B(M,N)\to N^{\oplus n}\to T:=\Img(\phi^\vee\otimes_BN)\to 0 \\
\Sigma' \quad:\quad &
0\to T\to N^{\oplus m}\to T':=\Coker(\phi^\vee\otimes_BN)\to 0
\end{aligned}
$$
where $\phi^\vee:B^{\oplus n}\to B^{\oplus m}$ is the transpose of
$\phi$. Notice that $T$ and $T'$ are $C$-modules of finite type.
By invoking again \cite[Ch.IV, Lemma 8.9.4.1]{EGAIV-3}, we
find $f\in A\setminus\{0\}$ such that $T_f$ and $T'_f$ are
flat $A_f$-modules; after replacing again $A,B,C,M,N$ and $\phi$
by their respective localizations, we may therefore assume that
$T$ and $T'$ are flat $A$-modules. In this case, the induced
sequences $A'\otimes_A\Sigma$ and $A'\otimes_A\Sigma'$ are
again short exact, for every $A$-algebra $A'$, whence a left
exact sequence
$$
0\to A'\otimes_A\Hom_B(M,N)\to
N'^{\oplus n}\xrightarrow{\phi'^\vee\otimes_{B'}N'}
N'^{\oplus m}
\qquad
\text{with $\phi':=\phi\otimes_BB'$}.
$$
But we have a natural identification :
$\Ker(\phi'^\vee\otimes_{B'}N')\isom\Hom_{B'}(M',N')$, and
the resulting isomorphism
$A'\otimes_A\Hom_B(M,N)\isom\Hom_{B'}(M',N')$ is the natural
map of (i).

(ii): Arguing as in the proof of (i), we reduce easily to the
case where $B$ and $C$ are finitely presented $A$-algebras,
and $M$ (resp. $N$) is a finitely presented $B$-module (resp.
$C$-module). Then we may find a $\Z$-subalgebra of finite
type $A_0\subset A$, an $A_0$-algebra of finite type $B_0$,
a $B_0$-algebra of finite type $C_0$, a $B_0$-module $M_0$ of
finite type, a $C_0$-module of finite type $N_0$, and a
$B_0$-linear map $\phi_0:M_0\to N_0$ with isomorphisms
$$
A\otimes_{A_0}B_0\isom B
\qquad
A\otimes_{A_0}C_0\isom C
\qquad
B\otimes_{B_0}M_0\isom M
\qquad
C\otimes_{C_0}N_0\isom N.
$$
which identify $B\otimes_{B_0}\phi_0$ with $\phi$. It suffices
then to exhibit $f\in A_0$ such that $(\Coker\,\phi_0)_f$ is a
free $A_{0,f}$-module; the latter follows easily from
\cite[Th.24.1]{Mat}.

(iii): If $x\in\Img\,(K\otimes_A\phi)$, we have $x\in\Img\,\phi_f$
for some $f\in A\setminus\{0\}$, and then
$x\otimes 1\in\Img(\kappa(\fp)\otimes_A\phi)$ for every
$\fp\in\Spec\,A_f$. Conversely, suppose that $U\subset\Spec\,A$ is a
dense subset such that $x\otimes 1\in\Img(\kappa(\fp)\otimes_A\phi)$
for every $\fp\in U$. We pick $f\in A\setminus\{0\}$ as in (ii), so
that $(\Coker\,\phi)_f$ is a free $A$-module. Clearly
$U\cap\Spec\,A_f$ is still dense in $\Spec\,A_f$, hence we may
replace $A,B,C,M,N$ and $\phi$ by $A_f,B_f,C_f,M_f,N_f$ and
$\phi_f$, and assume that $\Coker\,\phi$ is a free $A$-module,
say with basis $(e_i~|~i\in I)$. Suppose now that
$x\notin\Img\,(K\otimes_A\phi)$, so that the image
$\bar x\in\Coker\,\phi$ of $x$ does not vanish. Then there exists
a finite subset $I_0\subset I$ and a system $(f_i~|~i\in I_0)$ of
non-zero elements of $A$, with $\bar x=\sum_{i\in I_0}f_ie_i$. Pick
$j\in I_0$, and any $\fp\in U\cap\Spec\,A_{f_j}$; then
$\Coker\,(\kappa(\fp)\otimes_A\phi)$ is a free $\kappa(\fp)$-vector
space with basis $(1\otimes e_i~|~i\in I)$, and
$1\otimes x=\sum_{i\in I_0}\bar f_i\otimes e_i$, where
$\bar f_i\in\kappa(\fp)$ denotes the image of $f_i$, for every
$i\in I$. Especially, $\bar f_j\neq 0$, whence $1\otimes x\neq 0$,
a contradiction.
\end{proof}

\sset\subsubsection{}\label{subsec_recovered-from-R6}
Let $(P,+,0)$ be a fine and saturated monoid such that $P^\gp$
is torsion-free, and fix a Banach norm $\Vert\cdot\Vert$ on the
finite dimensional $\R$-vector space $P^\gp_\R:=\R\otimes_\Z P^\gp$.
Let also
$$
P^\gp_\R(\rho):=\{\gamma\in P^\gp_\R~|~\Vert\gamma\Vert\leq\rho\}
\qquad
\text{for all $\rho\in\R_+$}.
$$
Define also $P^{(\infty)}$ and $P_\Q$ as in
\eqref{subsec_intro-log-rings}, as well as the polyhedral cone
$P_\R:=\R_+\otimes_\Z P\subset P^\gp_\R$. Moreover, for every
$m\in\N\setminus\{0\}$, let
$$
P_{(m)}:=\{\gamma\in P_\Q~|~m\gamma\in P^{(\infty)}\}.
$$

\begin{lemma}\label{cl_epsilon-m}
With the notation of \eqref{subsec_recovered-from-R6}, there exists
an integer $m>0$ with $(p,m)=1$ such that the following holds. For
every $\delta>0$ there exists a real number $\eps>0$ such that for
every $\beta_1,\beta_2\in P^\gp_\R$ with
$\Vert\beta_1+\beta_2\Vert<\eps$ we may find
$\beta'_1,\beta'_2\in P_{(m)}^\gp$ with :
\begin{enumerate}
\item
$\beta'_1+\beta'_2=0$
\item
$\Vert\beta'_i-\beta_i\Vert<\delta$ for $i=1,2$
\item
$(P_\R-\beta_i)\cap P^\gp\subset(P_\R-\beta'_i)\cap P^\gp$ for $i=1,2$.
\end{enumerate}
\end{lemma}
\begin{proof} Since $P^{(\infty)\gp}$ is $p$-divisible,
$P_{(m)}^\gp=P_{(pm)}^\gp$ for every integer $m>0$, so we may always
arrange that $(p,m)=1$, if all the other conditions are already
fulfilled. Moreover, let $\rho>0$ such that every class of
$P^\gp_\R/P^\gp$ admits a representative $\beta\in P^\gp_\R$ with
$\Vert\beta\Vert\leq\rho$. Then it is easily seen that, after
replacing $\beta_i$ by $\beta_i+(-1)^i\gamma$ for $i=1,2$, with a
suitable $\gamma\in P^\gp$, we may assume that
$\beta_1\in P^\gp_\R(\rho)$. Now, for every $\gamma\in P^\gp_\R$,
let $\bar\Omega(\gamma)$ be the topological closure of the subset
$\Omega(P_\R,P^\gp\cap(P_\R-\gamma))$ in $P^\gp_\R$ (notation of
\eqref{subsec_for-almost-lengths}); in view of condition (i)
and of proposition \ref{prop_linear-constr-part}(v), condition
(iii) is implied by :
\set\begin{equation}\label{eq_replace-iii}
\beta'_1\in\bar\Omega(\beta_1)\cap(-\bar\Omega(\beta_2)).
\end{equation}
Moreover, suppose that
\set\begin{equation}\label{eq_replace-ii}
\Vert\beta'_1-\beta_1\Vert\leq\delta/2.
\end{equation}
Then
$$
\Vert\beta'_2-\beta_2\Vert=\Vert\beta'_2+\beta_1-\beta_1-\beta_2\Vert
\leq\Vert-\beta'_1+\beta_1\Vert+\Vert\beta_1+\beta_2\Vert<
\delta/2+\eps.
$$
Hence condition (ii) will hold, provided $\eps\leq\delta/2$.
Thus, we have to exhibit $m>0$ and $\eps>0$ such that, for every
$\beta_1,\beta_2\in P^\gp_\R$ with $\Vert\beta_1\Vert\leq\rho$
and $\Vert\beta_1+\beta_2\Vert<\eps$, there exists
$\beta'_1\in P_{(m)}^\gp$ fulfilling \eqref{eq_replace-iii} and
\eqref{eq_replace-ii}. Then, by proposition
\ref{prop_linear-constr-part}(i) and lemma \ref{lem_obama}, it
suffices to find $\beta'_1\in P^\gp_\R$ fulfilling these two
latter conditions. By way of contradiction, suppose that such
$\beta'_1$ cannot always be found : this means that there
exists a sequence
$\underline\beta:=((\beta_{1,k},\beta_{2,k})~|~k\in\N)$
of pairs of elements in $P^\gp_\R$, with
$\beta_{1,k}\in P^\gp_\R(\rho)$ for every $k\in\N$, and such that
\begin{enumerate}
\alphaenu
\item
$\Vert\beta_{1,k}+\beta_{2,k}\Vert<2^{-k}$ for every $k\in\N$
\item
$\bar\Omega(\beta_{1,k})\cap(-\bar\Omega(\beta_{2,k}))\cap
(P^\gp_\R(\delta/2)+\beta_{1,k})=\emptyset$ for every $k\in\N$.
\end{enumerate}
However, by proposition \ref{prop_linear-constr-part}(i),
after replacing $\underline\beta$ by a subsequence we may assume
that both $\bar\Omega(\beta_{1,k})$ and $\bar\Omega(\beta_{2,k})$
are independent of $k$. Since $P^\gp_\R(\rho)$ is a compact subset,
we may also assume that the sequence $(\beta_{1,k}~|~k\in\N)$
converges to an element $\beta'_1\in P^\gp_\R(\rho)$; then clearly
$(\beta_{2,k}~|~k\in\N)$ converges to $\beta'_2:=-\beta'_1$. Since
$\beta_{i,k}\in\bar\Omega(\beta_{i,k})$ for $i=1,2$ and every
$k\in\N$, we also have $\beta'_i\in\bar\Omega(\beta_{i,k})$ for
$i=1,2$. Lastly, we have $\Vert\beta'_1-\beta_{1,k}\Vert<\delta/2$
for every sufficiently large $k\in\N$; this contradicts (b), and the
claim follows.
\end{proof}

\begin{theorem}\label{th_diagonal-idempotent}
In the situation of \eqref{subsec_both-data}, suppose that
$(A,P,\beta)$ is regular. We have :

{\em(i)}\ \
$e_f\in\theta_{\underline f}(\sR(B,\Gamma_\Q))$.

{\em(ii)}\ \
Let $J\subset A$ be the radical ideal such that
$\Spec\,A/J=X\setminus U_f$. Then
$Je_f\subset\theta_{\underline f}(\sR(B,\Gamma))$.
\end{theorem}
\begin{proof} For every datum $\underline f:=(f,\beta,P)$
as in \eqref{subsec_both-data}, let
$\sX(\underline f):=\Gamma(U_f,\cB\otimes_{\cO_X}\cB)$; denote
also by $\sY(\underline f)\subset\sX(\underline f)$ the image
of $\theta_{\underline f}$, and set
$\sZ(\underline f):=Ae_f\subset\sX(\underline f)$. We notice :

\begin{claim}\label{cl_base-change-XYZ}
Let $A'$ be a noetherian ring, $g:A\to A'$ a flat ring
homomorphism; set $X':=\Spec\,A'$, $f':=A'\otimes_Af$ and 
$\underline f':=(f',P,g\circ\beta)$. Then $U_{f'}=X'\times_XU_f$,
and we have a natural isomorphism
$$
A'\otimes_A\sX(\underline f)\isom\sX(\underline f')
$$
that identifies $\sY(\underline f')$ and $\sZ(\underline f')$
respectively with $A'\otimes_A\sY(\underline f)$ and
$A'\otimes_A\sZ(\underline f)$.
\end{claim}
\begin{pfclaim} The first assertion follows from
\cite[Ch.IV, Prop.17.7.1]{EGA4}. Next, set $B':=A'\otimes_AB$;
since $B$ is an $A$-module of finite type, $B'$ is an $A'$-module
of finite presentation. Hence, the natural map
$A'\otimes_AB^\vee\to B'^\vee:=\Hom_{A'}(B',A')$ is an isomorphism
(details left to the reader). Then the second assertion follows
by a direct inspection of the constructions.
\end{pfclaim}

Now, assertion (i) comes down to the inclusion
$\sZ(\underline f)\subset\sY(\underline f)$, and it suffices
to check that $\sZ(\underline f)_\fp\subset\sY(\underline f)_\fp$
in $\sX(\underline f)_\fp$, for every $\fp\in\Spec\,A$. Let
$j_{(\fp)}:A\to A_\fp$ be the localization map, set
$f_\fp:=A_\fp\otimes_Af:A_\fp\to B_\fp$, and
$\underline f{}_\fp:=(f_\fp,P,\beta\circ j_{(\fp)})$; by claim
\ref{cl_base-change-XYZ}, we are then reduced to checking
that $\sZ(\underline f{}_\fp)\subset\sY(\underline f{}_\fp)$ in
$\sX(\underline f{}_\fp)$, for every such $\fp$. Thus, in order
to show (i), {\em we may assume that $(A,\fm)$ is a local
noetherian ring}. Likewise we may reduce assertion (ii) to
the case where $(A,\fm)$ is local. Next, let $A^\wedge$ be the
$\fm$-adic completion of $A$; since the completion map
$j:A\to A^\wedge$ is faithfully flat, we are reduced to checking
that $A^\wedge\otimes_A\sZ(\underline f)\subset
A^\wedge\otimes_A\sY(\underline f)$ in
$A^\wedge\otimes_A\sX(\underline f)$. Moreover, the log structure
attached to the chart $j\circ\beta:P\to A^\wedge$ is still
regular (theorems \ref{th_charact-log-regular} and
\ref{th_reg-generizes}); after invoking again claim
\ref{cl_base-change-XYZ}, we may thus {\em assume that $(A,\fm)$
is a complete noetherian local ring} in order to prove (i), and
a similar argument reduces the proof of (ii) as well to the
complete case.

Next, let $\fp:=\beta^{-1}(A^\times)$; then $\beta$ extends
uniquely to a morphism of monoids $\beta':P_\fp\to A$ (notation
of remark \ref{rem_localize-monoids}(i)); moreover, the log
structures on the Zariski site of $X$ induced by $\beta$ and
by $\beta'$ are naturally isomorphic, so the datum
$\underline f':=(f,P_\fp,\beta')$ still fulfills the conditions of
\eqref{subsec_both-data} and $(A,P_\fp,\beta')$ is regular.
Moreover, $(P_\fp)^\gp=P^\gp$ and $(P_\fp)_\Q^\gp=P^\gp_\Q$, hence
$(P_\fp)_\Q$ is still naturally $\Gamma$-graded. Furthermore,
for every $\gamma\in P^\gp_\Q$ we have
$$
\gr_{[\gamma]}((P_\fp)_\Q)=(\gr_{[\gamma]}(P_\Q))_\fp
\qquad\text{whence :}\qquad
A\otimes_{\Z[P_\fp]}\Z[\gr_{[\gamma]}((P_\fp)_\Q)]=
A_{[\gamma]}.
$$
Hence, the $A$-modules
$\sX(\underline f'),\sY(\underline f'),\sZ(\underline f')$ are
naturally identified with
$\sX(\underline f),\sY(\underline f),\sZ(\underline f)$.
Thus, in order to prove (i), it suffices to show that
$\sZ(\underline f')\subset\sY(\underline f')$; after replacing
$\underline f$ by $\underline f'$, we may therefore {\em assume
that the chart $\beta$ is local at the closed point of $X$} (see
definition \ref{def_chart}(vi)). Arguing likewise we reduce
as well the proof of (ii) to the case where $\beta$ is local.

Recall now that $P$ admits a decomposition
$P\isom Q\times P^\times$, where $Q$ is a fine, sharp and saturated
monoid, and $P^\times\subset P$ is the abelian group of invertible
elements of $P$ (lemma \ref{lem_decomp-sats}). Let $\alpha:Q\to A$
be the restriction of $\beta$; by direct inspection, we see that
the inclusion $Q\to P$ induces an isomorphism between the log
structures on $\Spec\,A$ associated with the charts $(P,\beta)$
and $(Q,\alpha)$ : see \eqref{subsec_up_and_down-log}; hence, the
datum $\underline f'':=(f,Q,\alpha)$ fulfills again the conditions
of \eqref{sec_diagonal-idempotents}, and $(A,Q,\alpha)$ is regular.
In light of lemma \ref{lem_finite-Sigma}(ii), we may then replace
$\underline f$ by $\underline f''$, and thus {\em assume that
$\beta$ is a sharp chart at the closed point of $X$} (see definition
\ref{def_chart}(vi)).

After these preliminaries, we may then assume, as in
\eqref{subsec_jeremy}, that there exists $r\in\N$, a coefficient
ring $(V,\fm_V)$ of $A$, and a surjective ring homomorphism
$$
\phi:V[[P\times\N^{\oplus r}]]\to A
$$
and $\beta$ is the composition of $\phi$ with the natural
inclusion map $P\to V[[P\times\N^{\oplus r}]]$.

\begin{claim}\label{cl_ok-if-A-is-F_p-alg}
The theorem holds if $A$ contains a field of positive characteristic.
\end{claim}
\begin{pfclaim} In this case, $V$ is the residue field of $A$,
and $\phi$ is an isomorphism. The induced morphism of monoids
$\beta':P\times\N^{\oplus r}\to V[[P\times\N^{\oplus r}]]\to A$ is
the chart for another regular log structure on the Zariski site
of $X$ (theorems \ref{th_charact-log-regular} and
\ref{th_reg-generizes}); in view of lemma \ref{lem_finite-Sigma}(ii)
we may then replace $\beta$ by $\beta'$, and assume that $A=V[[P]]$
and $\fm=\fm_P\cdot A$ (with $\fm_P:=P\setminus\{1\}$).
Next, let $k$ be an algebraic closure of $V$; the induced map
$V[P]\to k[P]$ is faithfully flat, so the same holds for its
completion $A\to k[[P]]$ (\cite[Th.22.4]{Mat}); moreover, the
inclusion map $P\to k[[P]]$ is again a chart for a regular log
structure on the Zariski site of $\Spec\,k[[P]]$. In view of
claim \ref{cl_base-change-XYZ}, we may then replace $V$ by $k$
and $B$ by $k[[P]]\otimes_AB$, and assume that $A=k[[P]]$. In
this case, notice that, with the notation of
\eqref{subsec_perfectification}, we have :
$$
A^\perf=\Gamma\times_{\Gamma_\Q}A_\infty.
$$
Let now $h\in J$. The isomorphism \eqref{eq_perfect-etale} yields
a natural identification of $A^\perf$-algebras :
$$
\Gamma\times_{\Gamma_\Q}B_\infty[h^{-1}]\isom B^\perf[h^{-1}]
$$
and then remark \ref{rem_perfectifcation} shows that $h^{1/p^n}e_f$ lies
in the image of the induced map
\set\begin{equation}\label{eq_perf-3-times}
B^\perf\otimes_{A^\perf}B^\perf=(B\otimes_AB)^\perf\to
B_\infty\otimes_{A_\infty}B_\infty[h^{-1}]\to
\Gamma(U_f,\cB_\infty\otimes_{\cO_{X_\infty}}\cB_\infty)
\end{equation}
for every $n\in\N$. On the other hand, $B^\perf$ is an integral
$A^\perf$-algebra, and $A^\perf$ is a normal domain, hence the
trace map of the finite \'etale map $A^\perf[h^{-1}]\to B^\perf[h^{-1}]$
restricts to an $A^\perf$-linear map $B^\perf\to A^\perf$ (lemma
\ref{lem_axiomatize-familiar}(i)) which in turns induces an
$A^\perf$-linear map
$$
\omega_{B^\perf/A^\perf}:
B^\perf\to\Hom_{A^\perf}(B^\perf,A^\perf)\to B^\vee_\infty
$$
such that
$$
A_\infty[h^{-1}]\otimes_{A^\perf}\omega_{B^\perf/A^\perf}=
\Gamma(\Spec\,A[h^{-1}],\omega_{\cB_\infty/\cO_{X_\infty}}).
$$
A simple inspection then shows that \eqref{eq_perf-3-times}
is the composition of
$\omega_{B^\perf/A^\perf}\otimes_{A^\perf}\omega_{B^\perf/A^\perf}$ with
\eqref{eq_was-theta_f}. Summing up, this proves that $h^{1/p^n}e_f$
lies in the image of the map \eqref{eq_was-theta_f} for every
$n\in\N$, and especially, $he_f\in\theta_{\underline f}(\sR(B,\Gamma))$,
whence assertion (ii).

In order to show (i), say that $h=\sum_{\gamma\in P}h_\gamma\cdot\gamma$
for a system $(h_\gamma~|~\gamma\in P)$ of elements of $k$, and
pick $\gamma_0\in P$ with $h_{\gamma_0}\neq 0$. According to claim
\ref{cl_good-gammas}, there exists $n\in\N$ such that
$\gr_{[\gamma_0^{1/p^n}]}P_\Q=\gamma_0^{1/p^n}P$; then we have an $A$-linear
isomorphism :
\set\begin{equation}\label{eq_legal-eagle}
\sX(\underline f)\isom\gr_{[\gamma_0^{1/p^n}]}
\Gamma(U_f,\cB_\infty\otimes_{\cO_{X_\infty}}\cB_\infty)
=\sX(\underline f)\otimes_AA_{[\gamma_0^{1/p^n}]}
\qquad
x\mapsto x\otimes\gamma_0^{1/p^n}.
\end{equation}
We have already remarked that $h^{1/p^n}e_f$ lies in the image
of \eqref{eq_was-theta_f}; let us write
$h^{1/p^n}=\sum_{[\gamma]\in\Gamma}h_{[\gamma]}$, with
$h_{[\gamma]}\in A_{[\gamma]}$ for every $[\gamma]\in\Gamma$; then
$h^{1/p^n}=\sum_{\gamma\in P}h^{1/p^n}_\gamma\cdot\gamma^{1/p^n}$, and
$h_{[\gamma_0^{1/p^n}]}=h_{\gamma_0}^{1/p^n}\cdot\gamma_0^{1/p^n}\cdot a$
for some $a\in 1+\fm$. Hence, $h_{[\gamma_0^{1/p^n}]}e_f$ is in the
image of
$\gr_{[\gamma_0^{1/p^n}]}(B^\vee_\infty\otimes_{A_\infty}B^\vee_\infty)$, and
so the same holds for $e_f\otimes\gamma_0^{1/p^n}$. Thus, there
exists a finite subset $\Sigma\subset\Gamma$ and an element of
$\bigoplus_{[\gamma]\in\Sigma}
B^\vee_{[\gamma]}\otimes_AB^\vee_{[\gamma_0^{1/p^n}/\gamma]}$ whose image
under the map \eqref{eq_was-theta_f} equals $e_f\otimes\gamma_0^{1/p^n}$.
Then, lemma \ref{cl_epsilon-m} yields for every
$[\gamma]\in\Sigma$ an element $[\gamma']\in\Sigma_\Q$ such that :
$$
\gamma^{-1}P_\Q\cap P^\gp\subset\gamma^{\prime-1}P_\Q\cap P^\gp
\qquad\text{and}\qquad
(\gamma/\gamma_0^{1/p^n})P_\Q\cap P^\gp\subset
\gamma'P_\Q\cap P^\gp.
$$
With $\Sigma':=\{[\gamma']~|~[\gamma]\in\Sigma\}$, we finally
deduce a commutative diagram of $A$-modules :
$$
\xymatrix{ 
\bigoplus_{[\gamma]\in\Sigma}
B^\vee_{[\gamma]}\otimes_AB^\vee_{[\gamma_0^{1/p^n}/\gamma]} \ar[r]  \ar[d] &
\sX(\underline f)\otimes_AA_{[\gamma_0^{1/p^n}]} \\
\sR(B,\Sigma')  \ar[r] &
\sX(\underline f) \ar[u]
}$$
whose right vertical arrow is the isomorphism \eqref{eq_legal-eagle}
and whose bottom horizontal arrow is the restriction of
$\theta_{\underline f}$; it follows easily that
$e_f\in\theta_{\underline f}(\sR(B,\Sigma'))$, as required.
\end{pfclaim}

{\em Suppose next that $A$ is a $\Q$-algebra}, so that $V$ is a
field of characteristic zero, and $\phi$ is an isomorphism. We pick
a complete discrete valuation ring $W$ with residue field $V$,
and a surjective map of monoids $\pi:\N^{\oplus r'}\to P$; for
$s:=r+r'$ we get a surjective map of $W$-algebras
$W[\N^{\oplus s}]\to A_0:=V[P\times\N^{\oplus r}]$ whose completion
is a surjective map $W[[\N^{\oplus s}]]\to A$. According to
\cite[\S3.6, Th.12]{Bo-Ray}, the completion map
$W[\N^{\oplus s}]\to W[[\N^{\oplus s}]]$ is the colimit of a filtered
system $(R_\lambda~|~\lambda\in\Lambda)$ of smooth
$W[\N^{\oplus s}]$-algebras; set
$A_\lambda:=R_\lambda\otimes_{W[\N^{\oplus s}]}A_0$ for every
$\lambda\in\Lambda$. Then $A$ is the colimit of the filtered
system of smooth $A_0$-algebras
$A_\bullet:=(A_\lambda~|~\lambda\in\Lambda)$. For every
$\lambda\in\Lambda$, let $X'_\lambda$ be the connected component
of $\Spec\,A_\lambda$ containing the image of $X$, and let
$A'_\lambda$ be the quotient of $A_\lambda$ with
$\Spec\,A'_\lambda=X'_\lambda$. It is easily seen that the colimit
of the induced system $A'_\bullet:=(A'_\lambda~|~\lambda\in\Lambda)$
is still $A$; after replacing $A_\bullet$ by $A'_\bullet$, we may
then assume that $A_\lambda$ is a domain for every $\lambda\in\Lambda$. 
Next, after replacing $\Lambda$ by a cofinal subset, we may assume
that for every $\lambda\in\Lambda$ there exists an open subset
$U_\lambda\subset X_\lambda$ such that $U_f=X\times_{X_\lambda}U_\lambda$,
and $U_\mu=X_\mu\times_{X_\lambda}U_\lambda$ for every $\mu\in\Lambda$
with $\mu\geq\lambda$. We may then find $\lambda\in\Lambda$ and a
finite morphism of schemes $\phi_\lambda:Y_\lambda\to X_\lambda$ with
an isomorphism of $X$-schemes $X\times_{X_\lambda}Y_\lambda\isom\Spec\,B$.
Then there exists $\mu\in\Lambda$ with $\mu\geq\lambda$ such that
$U_\mu\times_{X_\lambda}Y_\lambda$ is a finite \'etale $U_\mu$-scheme; set
$Y_\mu:=X_\mu\times_{X_\lambda}Y_\lambda$, and $B_\mu:=\cO_{Y_\mu}(Y_\mu)$.
By construction, the morphism $Y_\mu\to X_\mu$ is surjective; since
$A_\mu$ is a domain, it follows easily that the induced ring
homomorphism $f_\mu:A_\mu\to B_\mu$ is injective. Let
$\beta_\mu:P\to A_0\to A_\mu$ be the natural map; by example
\ref{ex_regular-log-ring}, the datum $(A_\mu,P,\beta_\mu)$ is a
regular log ring, and $\underline f{}_\mu:=(f_\mu,P,\beta_\mu)$
is a datum as in \eqref{subsec_both-data}. In light of remark
\ref{rem_functoriality-of-theta_f}(i), it then suffices to prove
assertions (i) and (ii) for the map $\theta_{\underline f{}_\mu}$;
the latter is covered by the following :

\begin{claim}
The theorem holds if there exists a field $K$ of characteristic
zero and a smooth ring homomorphism $\psi:K[P]\to A$, such that
$\beta$ is the composition of $\psi$ and the inclusion map
$P\to K[P]$ (see example \ref{ex_regular-log-ring}).
\end{claim}
\begin{pfclaim} The field $K$ is the colimit of the filtered
system of its smooth $\Z$-subalgebras
$(V_\lambda~|~\lambda\in\Lambda)$. Then we may find
$\lambda\in\Lambda$, a smooth $V_\lambda[P]$-algebra $A_\lambda$,
and a finite, injective and generically finite ring homomorphism
$f_\lambda:A_\lambda\to B_\lambda$ with isomorphisms
$$
K\otimes_{V_\lambda}A_\lambda\isom A
\qquad
A\otimes_{A_\lambda}B_\lambda\isom B
$$
(details left to the reader). Let $\alpha:P\to V_\lambda[P]\to A_\lambda$
be the natural map; since $V_\lambda$ is a regular noetherian ring,
the datum $(A_\lambda,P,\alpha)$ is a regular log ring, by example
\ref{ex_regular-log-ring}, and in light of remark
\ref{rem_functoriality-of-theta_f}(i), it suffices to show
assertions (i) and (ii) for the map $\theta_{\underline f{}_\lambda}$
corresponding to the datum $\underline f{}_\lambda:=(f_\lambda,P,\alpha)$.
Now, for every $\fp\in\Spec\,V_\lambda$, every $V_\lambda$-module
$M$, let $M(\fp):=\kappa(\fp)\otimes_{A_\lambda}M$; also, for every
homomorphism of $V_\lambda$-modules $h:M\to N$ let likewise
$h_{(\fp)}:=\kappa(\fp)\otimes_{V_\lambda}h:M(\fp)\to N(\fp)$;
moreover, let $\alpha_{(\fp)}:P\to A_\lambda(\fp)$ be the composition
of $\alpha$ with the natural map $A_\lambda\to A_\lambda(\fp)$.
Let also $a\in A_\lambda\setminus\{0\}$ such that
$\Spec\,A_\lambda[a^{-1}]\subset U_{f_\lambda}$.
According to \cite[Ch.IV, Prop.9.5.3]{EGAIV-3}, the set
$Z:=\{\fp\in\Spec\,V_\lambda~|~
\text{$\Spec\,A_\lambda[a^{-1}](\fp)$ is dense in $\Spec\,A_\lambda(\fp)$}\}$
is {\em constructible} in $\Spec\,V_\lambda$; furthermore, $Z$
contains the generic point of $\Spec\,V_\lambda$, so it is
{\em dense} in $\Spec\,V_\lambda$. For every $\fp\in Z$, the datum
$\underline f{}_{\lambda,(\fp)}:=(f_{\lambda,(\fp)},P,\alpha_{(\fp)})$
fulfills therefore the conditions of
\eqref{sec_diagonal-idempotents}, and moreover
$(A_\lambda(\fp),P,\alpha_{(\fp)})$ is still a regular log ring.
Set $X_\lambda:=\Spec\,A_\lambda$, and denote by $\cB_\lambda$ the
quasi-coherent $\cO_{X_\lambda}$-algebra arising from $B_\lambda$.
By remark \ref{rem_functoriality-of-theta_f}(i), for $\Delta$
either equal to $\Gamma_\Q$ or to $\Gamma$ we have a commutative
diagram :
$$
\xymatrix@C+20pt{
\sR(B_\lambda,\Delta) \ar[r]^-{\theta'_\Delta} \ar[d] &
B_\lambda\otimes_{A_\lambda}B_\lambda[a^{-1}] \ar[d] \\
\sR(B_\lambda(\fp,\Delta)) \ar[r]^-{\theta_\Delta^{\prime(\fp)}} &
B_\lambda\otimes_{A_\lambda}B_\lambda(\fp)[a^{-1}]
}$$
whose right vertical arrow maps $e_{f_\lambda}$ to $e_{f_{\lambda,(\fp)}}$,
and where $\theta'_\Delta$ is the composition of the restriction of
$\theta_{\underline f{}_\lambda}$ to $\sR(B,\Delta)$ with the restriction
map $\Gamma(U_{f_\lambda},\cB_\lambda\otimes_{\cO_{X_\lambda}}\cB_\lambda)\to
B_\lambda\otimes_{A_\lambda}B_\lambda[a^{-1}]$, and likewise
$\theta_\Delta^{\prime(\fp)}$ is deduced from
$\theta_{\underline f{}_\lambda,(\fp)}$. Lemma \ref{lem_finite-Sigma}(i)
yields a finite subset $\Sigma\subset\Delta$ such that
$$
\theta_{\underline f{}_{\lambda,(\fp)}}(\sR(B(\fp),\Sigma))=
\theta_{\underline f{}_{\lambda,(\fp)}}(\sR(B(\fp),\Delta))
\qquad
\text{for every $\fp\in Z$}.
$$
On the other hand, in view of lemma \ref{lem_from-generic-flatness}(i),
after replacing $V_\lambda$ by a suitable localization $V_\lambda[t^{-1}]$
(with $t\neq 0$), and $Z$ by $Z\cap\Spec\,V[t^{-1}]$, we may assume that
the natural maps
$$
(B^\vee_\lambda)_{[\gamma]}(\fp)\to(B_\lambda(\fp))^\vee_{[\gamma]}
\qquad
(B^\vee_\lambda)_{[1/\gamma]}(\fp)\to(B_\lambda(\fp))^\vee_{[1/\gamma]}
$$
are isomorphisms, for every $\fp\in Z$ and every $[\gamma]\in\Sigma$.
Summing up, we may assume that {\em the image of
$\theta_\Delta^{\prime(\fp)}$ equals the image of
$\theta'_{\Delta,(\fp)}$, for every $\fp\in Z$}. Lastly, recall
that the set $Z'$ of all $\fp\in\Spec\,V_\lambda$ such
that $\kappa(\fp)$ is a finite field is {\em dense in the
constructible topology of\/ $\Spec\,V_\lambda$}, since $V_\lambda$
is a $\Z$-algebra of finite type (\cite[Ch.IV, Cor.10.4.6]{EGAIV-3}).
Hence, $Z\cap Z'$ is a dense subset of $\Spec\,V_\lambda$, and by
claim \ref{cl_ok-if-A-is-F_p-alg}, the image of
$\theta_{\Gamma_\Q}^{\prime(\fp)}$ (resp. of $\theta_\Gamma^{\prime(\fp)}$)
contains $e_{f_{\lambda,(\fp)}}=1\otimes e_{f_\lambda}$ (resp.
$be_{f_{\lambda,(\fp)}}=1\otimes be_{f_\lambda}$) for every $\fp\in Z\cap Z'$
(resp. and every $b\in J$). By lemma
\ref{lem_from-generic-flatness}(iii) it follows that the
image of $\theta'_{\Gamma_\Q}$ contains $e_f$ (resp. that the
image of $\theta'_\Gamma$ contains $be_f$ for every $b\in J$).
The claim is an immediate consequence.
\end{pfclaim}

It remains to deal with the case where $A$ does not contain a
field, hence $A=V[[P]]/(h)$, for a discrete valuation ring $(V,\fm_V)$
whose residue field $k$ has positive characteristic, such that
$\fm_V=pV$, and with some $h:=\sum_{\gamma\in P}h_\gamma\cdot\gamma$
such that $h_0\in pV\setminus p^2V$. Arguing as in
\eqref{subsec_jeremy}, and taking into account claim
\ref{cl_base-change-XYZ}, we may then further reduce to the case
where $k$ is algebraically closed. In this situation, let
$g\in J\setminus\{0\}$, and set
$$
A_\Gamma:=\Gamma\times_{\Gamma_\Q}A_\infty
\qquad\text{and}\qquad
A'_\Gamma:=A_\Gamma[X^{1/p^\infty}]/(X-g)=A_\Gamma[g^{1/p^\infty}].
$$
We get a basic setup $(A_\Gamma,\fn)$, where $\fn\subset A_\Gamma$ is
the ideal generated by $(g^\delta~|~\delta\in\N[1/p]\setminus\{0\})$.
Moreover, let $D$ be the $p$-integral closure of $A'_\Gamma$ in
$A'_\Gamma[1/p]$, and $D^\wedge$ the $p$-adic completion of $D$; also,
set $D^\wedge_1:=(D^\wedge)^a_*$, the ring of almost elements of the
$(A_\Gamma,\fn)^a$-algebra $(D^\wedge)^a$, and denote by $B'_1$ the
integral closure of $D^\wedge_1$ in $B':=D^\wedge\otimes_AB[1/g]$.
Pick $N\in\N$ and $e\in B\otimes_AB$ whose image in
$B\otimes_AB[1/g]$ equals $g^Ne_f$, and let
$e'\in B'_1\hat\otimes_{D^\wedge}B'_1$ be the image of $e$.

\begin{claim}\label{cl_cleaning-up}
For every $\delta\in\N[1/p]\setminus\{0\}$ we have
$g^\delta e'\in g^N(B'_1\hat\otimes_{D^\wedge}B'_1)$.
\end{claim}
\begin{pfclaim} For every $n\in\N$, the $(D/p^nD)^a$-algebra
$(B'_1/p^nB'_1)^a$ is \'etale of finite rank, and we let
$\eps_n\in(B'_1/p^nB'_1\otimes_DB'_1/p^nB')^a_*$ be the corresponding
diagonal idempotent. The system $(\eps_n~|~n\in\N)$ corresponds to
a unique element
$$
\eps\in(B'_1\hat\otimes_{D^\wedge}B'_1)^a_*
$$
where $B'_1\hat\otimes_{D^\wedge}B'_1$ denotes the $p$-adic
completion of $B'_1\otimes_{D^\wedge}B'_1$. Let
$\mu:B'_1\otimes_{D^\wedge}B'_1\to B'_1$ be the multiplication
law of $B'_1$, and $\hat\mu$ the $p$-adic completion of $\mu$;
by construction, we see that :
\set\begin{equation}\label{eq_no-gronda}
\hat\mu(\eps)=1
\qquad\text{and}\qquad
\eps\cdot(1\otimes x-x\otimes 1)=0
\quad\text{for every $x\in B'_1=(B'_1)^a_*$}.
\end{equation}
By corollary \ref{cor_perf-Abhyankar}(ii), the completion map
$B'_1\otimes_{D^\wedge}B'_1\to B'_1\hat\otimes_{D^\wedge}B'_1$
induces an isomorphism
$$
B\otimes_AB\otimes_AD^\wedge[1/g]\isom
B'_1\otimes_{D^\wedge}B'_1[1/g]\isom
B'_1\hat\otimes_{D^\wedge}B'_1[1/g]
$$
and it follows easily from \eqref{eq_no-gronda} that the image
of $e_f\otimes 1$ in $B'_1\hat\otimes_{D^\wedge}B'_1[1/g]$ agrees
with the image of $\eps$. But $g$ is a regular element of
$(B'_1\hat\otimes_{D^\wedge}B'_1)^a_*$, due to corollary
\ref{cor_perf-Abhyankar}(ii), hence the image of $e$ in
$(B'_1\hat\otimes_{D^\wedge}B'_1)^a_*$ agrees with $g^N\eps$, whence
the claim.
\end{pfclaim}

Recall that $D^\wedge_1$ is integrally closed in $D^\wedge[1/g]$
(theorem \ref{th_perf-Abhyankar}(iv)), and $A$ is a normal
domain; then the trace maps for the \'etale ring homomorphisms
$D^\wedge[1/g]\to B'$ and $A[1/g]\to B[1/g]$ restrict to a
$D^\wedge_1$-linear map and respectively an $A$-linear map
$$
\tr_{B'_1/D^\wedge_1}:B'_1\to D^\wedge_1
\qquad
\tr_{B/A}:B\to A
$$
(lemma \ref{lem_axiomatize-familiar}(i)) whence an $A$-linear
map and a $D^\wedge_1$-linear map
$$
\begin{aligned}
\tau&:B\to B^\vee:=\Hom_A(B,A) & &
b\mapsto(b'\mapsto\tr_{B/A}(bb')) \\
\tau'&:B'_1\to\Hom_A(B,D^\wedge_1) & &
b\mapsto(b'\mapsto\tr_{B'_1/D^\wedge_1}(bb')).
\end{aligned}
$$
Let also $i:A\to A_\Gamma$ and $i_D:D^\wedge\to D^\wedge_1$ be the
natural maps; since $i^a_D$ is an isomorphism of $A_\Gamma^a$-modules,
we get by adjunction a unique $A_\Gamma$-linear map
$j_D:\tilde\fn\otimes_{A_\Gamma}D^\wedge_1\to D^\wedge$ such that
$i_D\circ j_D=\mu\otimes_{A_\Gamma}D^\wedge_1$, where as usual
$\tilde\fn:=\fn\otimes_{A_\Gamma}\fn$, and $\mu:\tilde\fn\to A$
is the multiplication map : $a\otimes a'\mapsto aa'$. Then, for
every $\delta,\delta'\in\N[1/p]\setminus\{0\}$ we let
$j_{\delta+\delta'}:D^\wedge_1\to D^\wedge$ be the $A_\Gamma$-linear map
such that $j_{\delta+\delta'}(d):=j_D(g^\delta\otimes g^{\delta'}\otimes d)$
for every $d\in D^\wedge_1$; it is easily seen that this map
depends only on the sum $\delta+\delta'$. Hence :
$$
i_D\circ j_\delta=g^\delta\cdot\one_{D^\wedge_1}
\qquad\text{and}\qquad
j_{\delta+\delta'}=g^\delta\cdot j_{\delta'}
\qquad
\text{for every $\delta,\delta'\in\N[1/p]\setminus\{0\}$}
$$
and moreover for every such $\delta$ we get a commutative
diagram of $A_\Gamma$-modules :
$$
\xymatrix@C+20pt{
A_\Gamma \ar[r]^-{g^\delta\cdot\one_{A_\Gamma}} \ar[d] &
A_\Gamma \ar[d] \\
D^\wedge_1 \ar[r]^{j_\delta} & D^\wedge
}$$
whose vertical arrows are the structure maps of the
$A_\Gamma$-algebras $D^\wedge_1$ and $D^\wedge$. Set
$$
B^\vee_\Gamma:=\Hom_A(B,A_\Gamma)
\qquad
B'^\vee_1:=\Hom_A(B,D^\wedge_1)
\qquad
C:=\Hom_A(B,D^\wedge)
$$
and let $i_*:B^\vee\to B^\vee_\Gamma$ (resp. $j_{\delta*}:B'^\vee_1\to C$)
be the $A$-linear (resp. $A_\Gamma$-linear) map such that
$i_*(\phi):=i\circ\phi$ for every $A$-linear map $\phi:B\to A$
(resp. $j_{\delta*}(\phi):=j_\delta\circ\phi$ for every $A$-linear map
$\phi:B\to D^\wedge_1$). To ease notation, set as well
$$
B^\vee_{\Gamma,i}:=\Hom_A(B,A_\Gamma/p^iA_\Gamma)
\qquad\text{and}\qquad
C_i:=\Hom_A(B,D/p^iD)
\qquad
\text{for every $i\in\N$}.
$$
We then obtain a commutative diagram for every
$\delta\in\N[1/p]\setminus\{0\}$ and every $i\in\N$ :
$$
\xymatrix{
B\otimes_AB \ar[rr]^-{\tau\otimes_A\tau} \ar[d] & &
B^\vee\otimes_AB^\vee
\ar[rr]^-{g^\delta\cdot i_*\otimes_Ag^\delta\cdot i_*} \ar[d] & &
B^\vee_\Gamma\otimes_{A_\Gamma}B^\vee_\Gamma \ar[r] \ar[d] &
B^\vee_{\Gamma,i}\otimes_{A_\Gamma}B^\vee_{\Gamma,i}
\ar[d] \\
B'_1\hat\otimes_{D^\wedge}B'_1
\ar[rr]^-{\tau'\hat\otimes_{D^\wedge}\tau'} & &
B'^\vee_1\hat\otimes_{D^\wedge}B'^\vee_1
\ar[rr]^-{j_{\delta*}\hat\otimes_{D^\wedge}j_{\delta*}} & &
C\hat\otimes_{D^\wedge}C \ar[r] & C_i\otimes_DC_i.
}$$
For given $\delta$, pick $\delta',\delta''\in\N[1/p]\setminus\{0\}$
with $\delta=\delta'+\delta''$; then
$$
(j_{\delta*}\circ\tau')\hat\otimes_{D^\wedge}(j_{\delta*}\circ\tau')(e')=
(j_{\delta'*}\circ\tau')\hat\otimes_{D^\wedge}(j_{\delta'*}\circ\tau')
(g^{2\delta''}e').
$$
By claim \ref{cl_cleaning-up}, we deduce that the image of
$(j_{\delta*}\circ\tau')\hat\otimes_{D^\wedge}(j_{\delta*}\circ\tau')(e')$
in $C_i\otimes_DC_i$ is divisible by $g^N$. Recall now that $D$ is
a faithfully flat $A_\Gamma$-algebra (theorem \ref{th_Zimbabwe});
since $A$ is noetherian, and $B$ is a finite $A$-algebra, we
deduce a natural isomorphism
$$
B^\vee_{\Gamma,i}\otimes_{A_\Gamma}D\isom C_i
\qquad
\text{for every $i\in\N$}.
$$
Set $e'':=(i_*\circ\tau)\otimes_A(i_*\circ\tau)(e)$; it follows that
the image of $ge''$ in $B^\vee_{\Gamma,i}\otimes_{A_\Gamma}B^\vee_{\Gamma,i}$
is divisible by $g^N$, for every $i\in\N$. Next, for every $i\in\N$
and every subset $\Delta\subset\Gamma_\Q$, let :
$$
\sR(B,\Delta,i):=\bigoplus_{[\gamma]\in\Delta}
B^\vee_{[\gamma],i}\otimes_AB^\vee_{[1/\gamma],i}
\quad\text{with $B^\vee_{[\gamma],i}:=\Hom_A(B,A_{[\gamma]}/p^iA_{[\gamma]})$
for every $[\gamma]\in\Gamma_\Q$}
$$
and notice that we have an $A$-linear surjection $\sR(B,\Gamma,i)\to
\gr_0(B^\vee_{\Gamma,i}\otimes_{A_\Gamma}B^\vee_{\Gamma,i})$, for the natural
$\Gamma$-grading on $B^\vee_{\Gamma,i}\otimes_{A_\Gamma}B^\vee_{\Gamma,i}$.
Summing up, for every $i\in\N$ there exists $e_i\in\sR(B,\Gamma,i)$
such that the image of $ge''$ in
$\gr_0(B^\vee_{\Gamma,i}\otimes_{A_\Gamma}B^\vee_{\Gamma,i})$ equals the
image of $g^Ne_i$. On the other hand, as in remark
\ref{rem_functoriality-of-theta_f}(ii) we obtain for every
$i\in\N$ a commutative diagram of $A$-linear maps
$$
\xymatrix{ A/p^iA\otimes_A\sR(B,\Gamma) \ar[r] \ar[d] &
A/p^iA\otimes_A\gr_0(B^\vee_\Gamma\otimes_{A_\Gamma}B^\vee_\Gamma)
\ar[d] \ar[r] & A/p^iA\otimes_A(B\otimes_AB)^\vee \ar[d] \\
\sR(B,\Gamma,i) \ar[r] &
\gr_0(B^\vee_{\Gamma,i}\otimes_{A_\Gamma}B^\vee_{\Gamma,i}) \ar[r] &
\Hom_A(B\otimes_AB,A/p^iA)
}$$
whose left vertical arrow is induced by the projections
$A_{[\gamma]}\to A_{[\gamma]}/p^iA_{[\gamma]}$ for every
$[\gamma]\in\Gamma_\Q$, and such that the composition of the
top horizontal arrows is $A/p^iA\otimes_A\xi_{\underline f}$.
Denote by $\xi_{\underline f,i}$ the composition of the bottom
horizontal arrows; arguing as in the proof of lemma
\ref{lem_finite-Sigma}(i) we find a finite subset
$\Sigma\subset\Gamma$ such that
$\xi_{\underline f}(\sR(B,\Gamma))=\xi_{\underline f}(\sR(B,\Sigma))$
and $\xi_{\underline f,i}(\sR(B,\Gamma,i))=
\xi_{\underline f,i}(\sR(B,\Sigma,i))$ for every $i\in\N$.
We notice:

\begin{claim}\label{cl_Hom-and-AR}
Let $R$ be a noetherian ring, $M,N$ two $R$-modules
of finite type, and $I\subset R$ an ideal. Set
$H^i:=\Hom_R(M,N/I^iN)$ for every $i\in\N$. Then there
exists $c\in\N$ such that
$$
\Img(H^{i+c}\to H^i)=\Img(\Hom_R(M,N)\to H^i)
\qquad
\text{for every $i\in\N$}.
$$
\end{claim}
\begin{pfclaim} Pick a finite presentation
$R^{\oplus q}\xrightarrow{\phi}R^{\oplus p}\to M$ of the $R$-module
$M$; let $\phi^\vee:R^{\oplus p}\to R^{\oplus q}$ be the transpose of
$\phi$, set $\phi^\vee_N:=\phi^\vee\otimes_AN:N^{\oplus q}\to N^{\oplus p}$
and $N':=\Img\,\phi^\vee_N$; by the Artin-Rees lemma
(\cite[Th.8.5]{Mat}) there exists $c\in\N$ such that
$N'\cap I^{i+c}N^{\oplus q}=I^i(N'\cap I^cN^{\oplus q})$ for
every $i\in\N$. We consider then the induced commutative diagram :
$$
\xymatrix{ 0 \ar[r] & \Hom_R(M,N) \ar[r] \ar[d] &
N^{\oplus p} \ar[rr]^-{\phi^\vee_N} \ar[d] & & N^{\oplus q} \ar[d] \\
0 \ar[r] & H^{i+c} \ar[r] \ar[d] &
N^{\oplus p}/I^{i+c}N^{\oplus p} \ar[rr]^-{\phi^\vee_{N/I^{i+c}N}} \ar[d] & &
N^{\oplus q}/I^{i+c}N^{\oplus q} \ar[d] \\
0 \ar[r] & H^i \ar[r] & N^{\oplus p}/I^iN^{\oplus p}
\ar[rr]^-{\phi^\vee_{N/I^iN}} & & N^{\oplus q}/I^iN^{\oplus q}
}$$
whose horizontal rows are left exact sequences. Let then
$\bar h\in\Ker\,\phi^\vee_{N/I^{i+c}N}$, and pick a representative
$h\in N^{\oplus p}$ for the class $\bar h$. Then
$\phi^\vee_N(h)\in N'\cap I^{i+c}N^{\oplus q}$, so there exists
$k\in\N$ such that $\phi^\vee_N(h)=\sum_{j=1}^ka_j\phi^\vee_N(x_j)$
for certain $a_1,\dots,a_k\in I^i$ and $x_1,\dots,x_k\in N^{\oplus p}$.
Set $h':=h-\sum_{j=1}^ka_jx_j$. Then $h'\in\Ker\,\phi^\vee_N$, and
the image of $h'$ in $N^{\oplus p}/I^iN^{\oplus p}$ agrees with the
image of $\bar h$, whence the claim.
\end{pfclaim}

By claim \ref{cl_Hom-and-AR}, for every $i\in\N$ there exists
$e'_i\in\sR(B,\Sigma)$ such that $\xi_{\underline f,i}(e_i)$ agrees
with the image of
$1\otimes\xi_{\underline f}(e'_i)\in A/p^iA\otimes_A(B\otimes_AB)^\vee$.
Now, let $M\subset(B\otimes_AB)^\vee$ be the $A$-submodule generated
by the system $(\xi_{\underline f}(e'_i)~|~i\in\N)$; we conclude that
the image of $ge''$ in $(B\otimes_AB)^\vee$ lies in the submodule
$g^NM+p^i(B\otimes_AB)^\vee$, for every $i\in\N$; but then it lies
already in $g^NM$, since $g^NM$ is a closed subset for the $p$-adic
topology of the finitely generated $A$-module $(B\otimes_AB)^\vee$
(\cite[Th.8.10(i)]{Mat}). Lastly, this shows that the image of
$ge''$ in $\sX(\underline f)$ lies in
$g^N\cdot\theta_{\underline f}(\sR(B,\Gamma))$; but this image equals
$g^{N+1}e_f$, whence $ge_f\in\theta_{\underline f}(\sR(B,\Gamma))$,
which achieves the proof of (ii).

To show (i), let $e''_i$ be the image of $e''$ in
$B^\vee_{\Gamma,i}\otimes_{A_\Gamma}B^\vee_{\Gamma,i}$ for every $i\in\N$,
and set
$$
I_i:=\{a\in A_\Gamma~|~ae''_i\in
g^N(B^\vee_{\Gamma,i}\otimes_{A_\Gamma}B^\vee_{\Gamma,i})\}
\qquad\text{and}\qquad
I_{D,i}:=\{d\in D~|~de''_i\in g^N(C_i\otimes_DC_i)\}.
$$
As already observed in the foregoing, $g^\delta\in I_{D,i}$ for
every $\delta\in\N[1/p]\setminus\{0\}$; then, since $D$ is a
flat $A_\Gamma$-algebra, we have $I_{D,i}=I_iD$ for every $i\in\N$
(details left to the reader), so that $I^k_{D,i}=I^k_iD$ for
every $i,k\in\N$. Since $D$ is a faithfully flat $A_\Gamma$-algebra,
it follows that $I^k_i=I^k_{D,i}\cap A_\Gamma$ for every $i,k\in\N$
(\cite[Th.7.5(ii)]{Mat}), whence $g\in I^k_i$ for every $i,k\in\N$,
and notice as well that $I_i$ is a $\Gamma$-graded ideal of $A_\Gamma$.
Let $\lambda_1,\dots,\lambda_k\in P^\vee:=\Hom_\Mnd(P_\R,\R_+)$
such that $\R\lambda_1,\dots,\R\lambda_k$ are the extremal
rays of $P^\vee_\R$, and set $\lambda:=\lambda_1+\cdots+\lambda_k$;
by lemma \ref{lem_explicit-epsilon} there exists a sequence
$(\gamma_n~|~n\in\N)$ of elements of $P^{(\infty)}$ with
$\lim_{n\to+\infty}\lambda(\gamma_n)=0$ and such that
$\gr_{[\gamma_n]}I_i=A_{[\gamma_n]}$ for every $n\in\N$, and by
claim \ref{cl_good-gammas} we may assume that
$\gr_{[\gamma_n]}P^{(\infty)}=\gamma_nP$ for every $n\in\N$, so
that $A_{[\gamma_n]}=A\cdot(1\otimes\gamma_n)$, and therefore
$$
(1\otimes\gamma_n)\cdot e''_i\in g^N\cdot
\gr_{[\gamma_n]}(B^\vee_{\Gamma,i}\otimes_{A_\Gamma}B^\vee_{\Gamma,i}).
$$
On the other hand, lemma \ref{cl_epsilon-m} yields for every
$[\gamma]\in\Gamma$ an element $[\gamma']\in\Sigma_\Q$ such that :
$$
\gamma^{-1}P_\Q\cap P^\gp\subset\gamma^{\prime-1}P_\Q\cap P^\gp
\qquad\text{and}\qquad
(\gamma/\gamma_n)P_\Q\cap P^\gp\subset\gamma'P_\Q\cap P^\gp.
$$
We then get a commutative diagram of $A$-modules :
$$
\xymatrix{ \sR(B,\Gamma_\Q,i) \ar[r] &
\gr_0(B^\vee_{\infty,i}\otimes_{A_\Gamma}B^\vee_{\infty,i}) \ar[r] &
\Hom_A(B\otimes_AB,A/p^iA) \ar[d] \\
\bigoplus_{[\gamma]\in\Gamma}
B^\vee_{[\gamma],i}\otimes_AB^\vee_{[\gamma_n/\gamma],i} \ar[u] \ar[r] &
\gr_{[\gamma_n]}(B^\vee_{\Gamma,i}\otimes_{A_\Gamma}B^\vee_{\Gamma,i}) \ar[r] &
\Hom_A(B\otimes_AB,A_{[\gamma_n]}/p^iA_{[\gamma_n]}) 
}$$
whose right vertical arrow is an isomorphism induced by scalar
multiplication by $1\otimes\gamma_n$. Moreover, the composition
of the top horizontal arrows is $\xi_{\underline f,i}$. With this
notation, we have
$e''_i\in\gr_0(B^\vee_{\Gamma,i}\otimes_{A_\Gamma}B^\vee_{\Gamma,i})$,
and it follows easily that the image of $e''_i$ lies in
$g^N\cdot\xi_{\underline f,i}(\sR(B,\Gamma_\Q,i))$, for every
$i\in\N$. Arguing as in the proof of lemma
\ref{lem_finite-Sigma}(i), we then find a finite subset
$\Delta\subset\Gamma_\Q$ such that the image of $e''_i$ in
$\Hom_A(B\otimes_AB,A/p^iA)$ lies in the $A$-submodule
$g^N\cdot\xi_{\underline f,i}(\sR(B,\Delta,i))$, for every
$i\in\N$. By invoking again claim \ref{cl_Hom-and-AR},
we deduce that for every $i\in\N$, the image of $e''$
in $(B\otimes_AB)^\vee$ lies in the $A$-submodule
$g^N\cdot\xi_{\underline f}(\sR(B,\Delta))+p^i\cdot(B\otimes_AB)^\vee$.
Arguing as in the foregoing proof of (ii), we conclude that
the image of $e''$ in $(B\otimes_AB)^\vee$ lies already in
$g^N\cdot\xi_{\underline f}(\sR(B,\Delta))$, and finally
$e_f\in\theta_{\underline f}(\sR(B,\Gamma_\Q))$, as stated.
\end{proof}

\sset\subsubsection{}\label{subsec_combinatorics}
Keep the notation of \eqref{subsec_both-data}, and for every
$m\in\N$, let $\Gamma_m:=\{\gamma\in\Gamma_\Q~|~m\gamma=0\}$.
By theorem \ref{th_diagonal-idempotent}, if $(A,P,\beta)$ is
regular, there exists $m\in\N$ such that
$e_f\in\theta_{\underline f}(\sR(B,\Gamma_m))$. It turns out
that a suitable integer $m$ with this property can be exhibited
purely in terms of the combinatorics of the monoid $P$. Indeed,
let $P^\vee$ be the dual of $P$ (see \eqref{subsec_dual-of-mon});
since $P$ is fine, $P^\vee$ is fine, sharp and saturated
(proposition \ref{prop_reflex-dual}(i,ii)), and we let
$F_1,\dots,F_k$ be the one-dimensional faces of $P^\vee$.
Then $F_i\simeq\N$ (theorem \ref{th_structure-of-satu}(ii)),
and we denote by $\lambda_i$ the unique generator of
$F_i$, for every $i=1,\dots,k$. Set
$\Lambda:=\{\lambda_1,\dots,\lambda_k\}$, and denote
by $\cE$ the set of all subsets $\Sigma\subset\Lambda$ such
that $\{\lambda^\gp\otimes_\Z\Q:P^\gp_\Q\to\Q~|~\lambda\in\Sigma\}$
is a basis of $(P^\gp_\Q)^\vee$. For every $\Sigma\in\cE$, let
$L_\Sigma\subset P^{\vee\gp}$ be the subgroup generated by
$(\lambda^\gp~|~\lambda\in\Sigma)$, and denote by $m_\Sigma$
the exponent of the finite abelian group $P^{\vee\gp}/L_\Sigma$.
Lastly, denote by $m_P$ the least common multiple of the finite
system of integers $(m_\Sigma~|~\Sigma\in\cE)$. We may then state :

\begin{corollary}\label{cor_diff-abhyankhar}
In the situation of theorem {\em\ref{th_diagonal-idempotent}},
we have $e_f\in\theta_{\underline f}(\sR(B,\Gamma_{m_p}))$.
\end{corollary}
\begin{proof} For every $\Sigma\in\cE$, let
$\lambda_\Sigma:P^\gp_\Q\isom\Q^\Sigma$ be the $\Q$-linear isomorphism
given by the rule :
$$
\gamma\mapsto(\lambda^\gp_\Q(\gamma)~|~\lambda\in\Sigma)
\qquad\text{for every $\gamma\in P^\gp_\Q$}.
$$
Notice that $P^\gp\subset\lambda^{-1}_\Sigma(\Z^\Sigma)$ for
every such $B$, and $m_\Sigma$ equals the exponent of the
finite abelian group $\lambda_\Sigma^{-1}(\Z^\Sigma)/P^\gp$.
With this notation, we have :

\begin{claim}\label{cl_refine-estimate}
For every $\beta\in P^\gp_\Q$ there exist $\Sigma\in\cE$ and
an element $\beta^\dagger\in\lambda^{-1}_\Sigma(\Z^\Sigma)$ such that
\set\begin{equation}\label{eq_beta-gamma}
\beta P_\Q\cap P^\gp\subset\beta^\dagger P_\Q\cap P^\gp
\qquad\text{and}\qquad
\beta^{-1}P_\Q\cap P^\gp\subset\beta^{\dagger-1}P_\Q\cap P^\gp.
\end{equation}
\end{claim}
\begin{pfclaim} Notice that, for every $\beta\in P^\gp_\R$,
we have
$$
\beta P_\R\cap P^\gp=\{x\in P^\gp~|~
\lambda^\gp(x)\geq\lambda^\gp_\R(\beta)\
\text{for every $\lambda\in\Lambda$}\}.
$$
For every such $\beta$ and every $\lambda\in\Lambda$, denote
by $N_{\beta,\lambda}$ (resp. $N'_{\beta,\lambda}$) the largest (resp.
smallest) integer which is smaller than (resp. larger than) or
equal to $\lambda^\gp_\R(\beta)$. Recalling that
$\lambda^\gp(P^\gp)\subset\Z$ for every $\lambda\in\Lambda$, it
is easily seen that \eqref{eq_beta-gamma} holds if and only if
\set\begin{equation}\label{eq_equiv-to-beta-gamma}
\lambda^\gp_\Q(\beta^\dagger)\in[N_{\beta,\lambda},N'_{\beta,\lambda}]
\qquad
\text{for every $\lambda\in\Lambda$}.
\end{equation}
Set $\Sigma_\beta:=\{\lambda\in\Lambda~|~\lambda^\gp_\R(\beta)\in\Z\}$.
Now, if $\beta\in P^\gp_\Q$ and $\Sigma_\beta$ spans $(P^\gp_\Q)^\vee$,
obviously the claim holds with $\beta^\dagger:=\beta$. Thus, a simple
induction argument reduces to showing that, if $\beta\in P^\gp_\Q$
and $\Sigma_\beta$ does not span $(P^\gp_\Q)^\vee$, there exists
$\gamma\in P^\gp_\Q$ such that \eqref{eq_equiv-to-beta-gamma}
holds, and such that $\Sigma_\gamma$ strictly contains $\Sigma_\beta$.
To this aim, let $l:P^\gp_\R\to\R^\Lambda$ be the $\R$-linear
map such that
$$
l(\gamma):=(\lambda^\gp_\R(\gamma)~|~\lambda\in\Lambda)
\qquad
\text{for every $\gamma\in P^\gp_\R$}
$$
and set
$$
K:=\prod_{\lambda\in\Lambda}[N_{\beta,\gamma},N'_{\beta,\gamma}]
\subset\R^\Lambda
\qquad
M:=\bigcap_{\lambda\in\Sigma_\beta}\Ker\,\lambda^\gp_\Q\subset P^\gp_\Q.
$$
Notice that $l$ is a closed immersion, $K$ is a compact
subset of $\R^\Lambda$, and by assumption, $M\neq 0$.
Clearly $\beta\in l^{-1}K$, and we pick any non-zero $\mu\in M$.
It follows that the subset
$$
\{a\in\R~|~\beta+a\mu\in l^{-1}K\}
$$
is non-empty and compact, so it admits a largest element $a_0$.
Set $\gamma:=\beta+a_0\mu$; by construction,
$\Sigma_\beta\subset\Sigma_\gamma$, and it is easily seen that there
must exist some $\lambda\in\Lambda\setminus\Sigma_\beta$ such that
$\lambda^\gp_\R(\gamma)\in\{N_{\beta,\lambda},N'_{\beta,\lambda}\}$,
since otherwise we could find a real number $a>a_0$ with
$\beta+a\mu\in l^{-1}K$. Thus $\Sigma_\gamma$ strictly contains
$\Sigma_\beta$. Lastly, since both $\beta$ and $\mu$ lie in
$P^\gp_\Q$ and $\lambda^\gp_\R(\beta)\notin\Z$, we must have
$\lambda^\gp_\R(a_0\mu)\in\Q\setminus\Z$, and since
$\lambda^\gp_\R(\mu)\in\Q$, we conclude that $a_0\in\Q$, hence
$\gamma\in P^\gp_\Q$, whence the claim.
\end{pfclaim}

Now, for every $[\beta]\in\Gamma_\Q$ pick $\beta^\dagger\in P^\gp_\Q$
fulfilling the condition of claim \ref{cl_refine-estimate}, and set
$\Delta:=\{[\beta^\dagger]~|~[\beta]\in\Gamma_\Q\}$. Notice that
$\Delta\subset\Gamma_{m_P}$. Arguing as in the proof of lemma
\ref{lem_finite-Sigma}(ii) we easily deduce that
$\theta_{\underline f}(\sR(B,\Delta))=
\theta_{\underline f}(\sR(B,\Gamma_\Q))$, whence the corollary.
\end{proof}

\begin{example}
(i)\ \
Suppose that $\dim P=2$. Then we may find a basis
$(f_1,f_2)$ of $P^\gp$ and integers $a,b\in\N$ such
that $b>a$, $(a,b)=1$ and $P$ is the submonoid
of $P^\gp$ generated by $f_1$ and $af_1+bf_2$ (example
\ref{ex_satu-dim-two}(ii)). With this notation, a
simple calculation shows that $m_P=b$ (details left
to the reader).

(ii)\ \
In the situation of theorem \ref{th_diagonal-idempotent}, suppose
that $P=\N^{\oplus r}\times\Z^{\oplus s}$ for some $r,s\in\N$; by
corollary \ref{cor_more-precisely}, the latter holds if and only
if $A$ is a regular ring. In this case, the integer $m_P$ of
\eqref{subsec_combinatorics} equals $1$, and therefore corollary
\ref{cor_diff-abhyankhar} says that
$e_f\in\theta_{\underline f}(B^\vee\otimes_AB^\vee)$. However, the latter
follows also more directly from lemma \ref{lem_finite-Sigma}(ii).
Moreover, in case $A$ is regular, the submodule
$\theta_{\underline f}(B^\vee\otimes_AB^\vee)$ can be interpreted in
terms of the {\em different ideal} of the finite ring extension
$f:A\to B$. We conclude this section with a review of different
ideals in the more general context of quasi-coherent algebras on
schemes; especially, remark \ref{rem_Hulot-demissionne} explains
the connection with theorem \ref{th_diagonal-idempotent}.
\end{example}

\sset\subsubsection{}\label{subsec_inverse-diff}
Let $X$ be a reduced scheme, $\cA$ a quasi-coherent
and finite $\cO_{\!X}$-algebra and $j:U\to X$ a
quasi-compact open immersion with dense image, such
that $j^*\cA$ is a locally free $\cO_{\!U}$-module.
Denote by
$$
\Max\,X
$$
the set of maximal points of $X$, and notice that
$\Max\,X\subset U$, by proposition
\ref{prop_closed-under-spec}(i). In this situation, there
is a trace form as in remark \ref{rem_idemp-and-traces}(iii) :
$$
t_{j^*\cA}:j^*\cA\otimes_{\cO_{\!U}}j^*\cA\to\cO_{\!U}
$$
whence a pairing
$$
t_{\cA,U}:\cA\otimes_{\cO_{\!X}}j_*j^*\cA\to
j_*(j^*\cA\otimes_{\cO_{\!U}}j^*\cA)
\xrightarrow{\ j_*t_{j^*\cA}\ } j_*\cO_{\!U}
$$
which in turns yields a morphism of $\cO_{\!X}$-modules
$$
\tau_{\cA,U}:j_*j^*\cA\to\cHom_{\cO_{\!X}}(\cA,j_*\cO_{\!U}).
$$
Since $X$ is reduced, the natural morphism
$\cO_{\!X}\to j_*\cO_{\!U}$ is injective, and
therefore the same holds for the induced
$\cO_{\!X}$-linear morphism
$$
\cA^\vee\to\cHom_{\cO_{\!X}}(\cA,j_*\cO_{\!U})
$$
so we may define
$$
\cD^{-1}_{\cA,U}:=\tau^{-1}_{\cA,U}(\cA^\vee).
$$
We call $\cD^{-1}_{\cA,U}\subset j_*j^*\cA$ the {\em inverse
different of $\cA$ relative to the open subset $U$}.

\begin{remark}\label{rem_taus-and-Us}
(i)\ \
In the situation of \eqref{subsec_inverse-diff}, notice
that the natural isomorphism
$$
\cHom_{\cO_{\!X}}(\cA,j_*\cO_{\!U})\isom j_*(j^*\cA)^\vee
$$
(see \eqref{eq_std-adjs}) identifies $\tau_{\cA,U}$ with
$j_*\tau_{j^*\cA}$, where $\tau_{j^*\cA}:j^*\cA\to j^*\cA^\vee$
is the $\cO_{\!U}$-linear morphism deduced from the trace
pairing $t_{j^*\cA}$.

(ii)\ \
Moreover, set $\cB:=\Img\,(\cA\to j_*j^*\cA)$. Since the
natural map $\cO_{\!X}\to j_*\cO_{\!U}$ is a monomorphism,
it is easily seen that the the epimorphism $\cA\to\cB$
induces an isomorphism $\cB^\vee\isom\cA^\vee$, and taking
(i) into account, we deduce that the latter induces
an isomorphism
$$
\cD^{-1}_{\cA,U}\isom\cD^{-1}_{\cB,U}.
$$

(iii)\ \
Suppose next that $U'\subset U$ is another open
subset, such that the open immersion $j':U'\to X$
is also quasi-compact with dense image. A direct
inspection of the definitions yields a commutative
diagram
$$
\xymatrix{ j_*j^*\cA \ar[rr]^-{\tau_{\cA,U}} \ar[d] & &
\cHom_{\cO_{\!X}}(\cA,j_*\cO_{\!U}) \ar[d] \\
j'_*j'{}^*\cA \ar[rr]^-{\tau_{\cA,U'}} & &
\cHom_{\cO_{\!X}}(\cA,j'_*\cO_{\!U'})
}$$
whose vertical arrows are monomorphisms, under the
current assumptions. It follows easily that
$$
\cD^{-1}_{\cA,U}=\cD^{-1}_{\cA,U'}\cap j_*j^*\cA.
$$

(iv)\ \
In the situation of (iii), suppose that $\cA$ is generically
\'etale, and let $U'\subset U$ be the largest open subset such
that $\cA_{|U'}$ is an \'etale $\cO_{U'}$-algebra. Then the
inclusion $U'\to U$ is a dense and quasi-compact open immersion.
For the proof, we may assume that $U$ is affine, say $U=\Spec\,R$,
and that $j^*\cA$ is the quasi-coherent $\cO_U$-algebra associated
to an $R$-algebra $A$ which is a projective $R$-module of finite
rank; then we easily reduce to the case where $A$ is a free
$R$-module of finite rank. In this case, let $(e_i~|~i=1,\dots,r)$
be a basis of the $R$-module $A$; consider the matrix
$D:=(t_{j^*\cA}(e_i,e_j)~|~i,j=1,\dots,r)$, and set $d:=\det(D)$.
We have $U'=\Spec\,R[d^{-1}]$, a quasi-compact dense open subset
of $U$ (\cite[Th.4.1.14]{Ga-Ra}), as required.

(v)\ \
In the situation of (iii), suppose additionally that
$j^*\cA$ is an \'etale $\cO_{\!U}$-algebra. Then $t_{j^*\cA}$
and $t_{j'{}^*\cA}$ are perfect pairings (\cite[Th.4.1.14]{Ga-Ra}).
Taking (i) into account, we conclude that $\tau_{\cA,U}$
is an isomorphism, and the same applies to $\tau_{\cA,U'}$,
so the inclusion map $j_*j^*\cA\to j'_*j'{}^*\cA$ restricts
to a natural isomorphism
$$
\cD^{-1}_{\cA,U}\isom\cD^{-1}_{\cA,U'}.
$$
More generally, if $j:U\to X$ and $j':U'\to X$ are any two
quasi-compact dense open immersions such that $j^*\cA$ and
$j'^*\cA$ are respectively an \'etale $\cO_U$-algebra and
an \'etale $\cO_{U'}$-algebra, then both $\cD^{-1}_{\cA,U}$
and $\cD^{-1}_{\cA,U'}$ are naturally identified with
$\cD^{-1}_{\cA,U\cap U'}$, and in this sense we may say that
the inverse different of a finite and generically \'etale
$\cO_{\!X}$-algebra is independent of the open subset $U$. We
shall therefore henceforth denote just by $\cD^{-1}_\cA$ the
inverse different of $\cA$.
\end{remark}

\begin{lemma}\label{lem_shortened}
In the situation of \eqref{subsec_inverse-diff}, denote
by $\cA^\nu$ the normalization of $j^*\cA$ over $X$ (see
remark {\em\ref{rem_almost-pairs}(ii)}). The following holds :
\begin{enumerate}
\item
If $X$ is normal, $\cA\subset\cD^{-1}_{\cA,U}$.
\item
Suppose that $\cA$ is a generically \'etale
$\cO_{\!X}$-algebra. Then :
\begin{enumerate}
\item
There exists a dense quasi-compact open subset $U'\subset U$ such
that the map $\tau_{\cA,U'}$ restricts to an $\cO_{\!X}$-linear
isomorphism $\cD^{-1}_\cA\isom\cA^\vee$.
\item
If $X$ is locally coherent, $\cD^{-1}_\cA$ is a reflexive
$\cO_{\!X}$-module of finite type.
\item
If $X$ is normal, $\cA^\nu\subset\cD^{-1}_\cA$.
\end{enumerate}
\item
Suppose that $X$ is normal, and $j^*\cA$ is an \'etale
$\cO_U$-algebra. Let also $j':U'\to X$ be a quasi-compact
open immersion with dense image, with $U'\subset U$, and
$\cA'^\nu$ the integral closure of $\cO_X$ in $j'^*\cA$.
Then the natural map $\cA^\nu\to\cA'^\nu$ is an isomorphism.
\end{enumerate}
\end{lemma}
\begin{proof} (i) follows immediately from lemma
\ref{lem_axiomatize-familiar}(i).

(ii.a) follows immediately from the discussion of
remark \eqref{rem_taus-and-Us}(iv,v), and (ii.b) follows
from (ii.a) and lemma \ref{lem_dual-is-rflx}.

(ii.c): Taking into account remark \ref{rem_taus-and-Us}(ii),
we may assume without loss of generality that
$\cA=\Img\,(\cA\to j_*j^*\cA)$, in which case the natural map
$\cA\to\cA^\nu$ is a monomorphism. Moreover, the assertion is
clearly local on $X$, hence we may assume that $X$ is affine,
in which case $\cA^\nu$ is the union of the filtered family of
its $\cA$-subalgebras $\cB$ that are quasi-coherent and finite
$\cO_{\!X}$-algebras. For any such $\cB$, the $\cO_{\!X}$-module
$\cB/\cA$ is supported on $X\setminus U$. It follows that
$(\cB/\cA)^\vee=0$, so the transpose map $\cB^\vee\to\cA^\vee$
is a monomorphism as well; by the same token, we see that
$\tau_{\cA,U}=\tau_{\cB,U}$ (notation of \eqref{subsec_inverse-diff}),
so $\cD^{-1}_\cB\subset\cD^{-1}_\cA$, whence $\cB\subset\cD^{-1}_\cA$,
by virtue of (i). Since $\cB$ is arbitrary, the assertion
follows.

(iii): We may assume that $X$ is local, say $X=\Spec\,R$
for a local and normal domain $R$, and let $K$ be the field
of fractions of $R$; then $A^\nu:=\cA^\nu(X)$ (resp. $A'^\nu$) is
the integral closure of $R$ in $\cA(U)$ (resp. in $\cA(U')$);
but under the stated assumptions, both $A^\nu$ and $A'^\nu$
are also the integral closure of $R$ in $K\otimes_R\cA(X)$,
whence the assertion.
\end{proof}

\begin{remark} Let $X$ be a reduced and quasi-separated
scheme, such that $\Max\,X$ is finite.

(i)\ \
Let also $\cF$ be a quasi-coherent $\cO_{\!X}$-module
of finite type. Then we claim that there exists a
dense and quasi-compact open immersion $j:U\to X$
such that $j^*\cF$ is a locally free $\cO_{\!U}$-module.
Indeed, for any maximal point $\eta\in X$, pick an affine
open subset $V_\eta\subset X$ such that $\eta$ is the
unique maximal point of $V_\eta$; hence $V_\eta=\Spec\,A$
for an integral domain $A$ and $\cF_{|V}=M^\sim$ for an
$A$-module $M$ of finite type. Set $K:=\Frac\,A$; then
$P:=M\otimes_AK$ is a finite dimensional $K$-vector space,
and we may find a element $f\in A\!\setminus\!\{0\}$, a
free $A_f$-module $Q$ of finite rank and an $A_f$-linear map
$\phi:Q\to M\otimes_AA_f$ such that
$\Ker\,\phi\otimes_AK=0=\Coker\,\phi\otimes_AK$.
However, the natural map $Q\to Q\otimes_{A_f}K$ is
injective, so we see that $\Ker\,\phi=0$, and on the
other hand $\Coker\,\phi$ is an $A_f$-module of finite
type, so we may find $g\in A\!\setminus\!\{0\}$ such
that $\phi_g$ is surjective. Set $U_\eta:=\Spec\,A_{fg}$;
by construction $\cF_{|U_\eta}$ is a locally free
$\cO_{\!U_\eta}$-module, and since $\eta$ is arbitrary,
the assertion follows easily.

(ii)\ \
Let $\cA$ be a quasi-coherent and finite $\cO_{\!X}$-algebra,
and suppose that $\cA_\eta$ is an \'etale $\cO_{\!X,\eta}$-algebra,
for every maximal point $\eta$ of $X$. Then we claim that there
exists an dense quasi-compact open immersion $j:U\to X$ such
that $j^*\cA$ is an \'etale $\cO_{\!U}$-algebra. Indeed, by (i),
we may already assume that $\cA$ is a locally free
$\cO_{\!X}$-module, and then we may argue as in remark
\ref{rem_taus-and-Us}(iv).
\end{remark}

\sset\subsubsection{}\label{subsec_jumping}
Keep the situation of \eqref{subsec_inverse-diff},
let $X'$ be another reduced scheme, and $g:X'\to X$
a quasi-compact and quasi-separated morphism of schemes
which restricts to a map
$$
\Max\,X'\to\Max\,X.
$$
Set
$$
U':=g^{-1}U
\qquad
\cA':=g^*\cA
$$
and denote by $j':U'\to X'$ the resulting open immersion.
It is easily seen that $j'$ is quasi-compact and has
dense image, and clearly $j'{}^*\cA'$ is a locally free
$\cO_{\!U'}$-module, so there is a well defined inverse
different $\cD^{-1}_{\cA',U'}$ relative to $U'$. Moreover,
we have a commutative diagram
\set\begin{equation}\label{eq_g-differ-pback}
{\diagram
g^*j_*j^*\cA
\ar[rrr]^-{g^*\tau_{\cA,U}} \ar[d]_a & & &
g^*\cHom_{\cO_{\!X}}(\cA,j_*\cO_{\!U}) \ar[d]_b &
\ar[l] g^*(\cA^\vee) \ar[d]^c \\
j'_*j'{}^*\cA' \ar[rrr]^-{\tau_{\cA',U'}} & & &
\cHom_{\cO_{\!X'}}(\cA',j'_*\cO_{\!U'}) &
\ar[l] \cA^{\prime\vee}
\enddiagram}
\end{equation}
(\cite[Lemma 4.1.13(iii)]{Ga-Ra}). It follows
that the left vertical arrow restricts to an
$\cO_{\!X'}$-linear map
\set\begin{equation}\label{eq_diff-base-change}
g^*\cD^{-1}_{\cA,U}\to\cD^{-1}_{\cA',U'}.
\end{equation} 

\begin{lemma}\label{lem_jumping}
The map \eqref{eq_diff-base-change} is an isomorphism
in the following cases :
\begin{enumerate}
\alphaenu
\item
if $\cA$ is a locally free $\cO_{\!X}$-module of finite type
and $j^*\cA$ is an \'etale $\cO_U$-algebra
\item
or if $g$ is a flat morphism.
\end{enumerate}
\end{lemma}
\begin{proof} For case (a), notice that when $\cA$ is locally
free, the map $c$ is an isomorphism, and if $j^*\cA$ is an
\'etale $\cO_U$-algebra, the discussion of remark
\ref{rem_taus-and-Us}(v) identifies \eqref{eq_diff-base-change}
with the map $c$.

If $g$ is flat, the map $a$ in \eqref{eq_g-differ-pback}
is an isomorphism (corollary \ref{cor_base-change-where}).
The isomorphism of remark \ref{rem_taus-and-Us}(i) identifies
$b$ with the composition
$$
g^*j_*\cHom_{\cO_{\!X}}(j^*\cA,\cO_{\!U})\xrightarrow{\ d\ }
j'_*g^*\cHom_{\cO_{\!X}}(j^*\cA,\cO_{\!U})\xrightarrow{\ e\ }
j'_*\cHom_{\cO_{\!U'}}(j'{}^*\cA',\cO_{\!U'})
$$
where $d$ is an isomorphism, since $g$ is flat
(corollary \ref{cor_base-change-where}), and the
same holds for $e$, since $j^*\cA$ is locally free
of finite type (proposition \ref{prop_replace-prop.12.3.5}(ii)).
To conclude, it then suffices to show that $c$ is
an isomorphism as well. The assertion is local on
$X$, hence we may assume that $X$ is affine; in this
case, we may find a finitely presented $\cO_{\!X}$-module
$\cG$ and an $\cO_{\!X}$-linear epimorphism
$\phi:\cG\to\cA$ such that $j^*\phi$ is an isomorphism. Set
$\cG':=g^*\cG$; there follows a commutative diagram
$$
\xymatrix{
g^*(\cA^\vee) \ar[rr]^-{g^*(\phi^\vee)} \ar[d] & &
g^*(\cG^\vee) \ar[d] \\
\cA'{}^\vee \ar[rr]^-{(\phi^*g)^\vee} & & \cG'{}^\vee
}$$
whose right vertical arrows is an isomorphism, by
virtue of \cite[Lemma 2.4.29(i)]{Ga-Ra}. Moreover
$\Ker\,(\phi^\vee)=0$ and
$\Coker\,(\phi^\vee)\subset(\Ker\,\phi)^\vee$,
and since $\Max\,X\subset U$, we see that
$\Ker\,\phi_\eta=0$ for every $\eta\in\Max\,X$,
whence $(\Ker\,\phi)^\vee=0$, since $X$ is reduced.
Thus, the top horizontal arrow of the foregoing
diagram is also an isomorphism; the same argument
applies as well to the bottom horizontal arrow, and
the assertion follows.
\end{proof}

\sset\subsubsection{}\label{subsec_jump-again}
In the situation of \eqref{subsec_jumping}, suppose now that
$X$ and $X'$ are locally coherent and $\cA$ is generically
\'etale. Taking into account lemma \ref{lem_shortened}(ii.b),
we see that \eqref{eq_diff-base-change} factors as a
composition
$$
g^*\cD^{-1}_\cA\xrightarrow{\ \beta\ }
(g^*\cD^{-1}_\cA)^{\vee\vee}\xrightarrow{\ \gamma\ }
\cD^{-1}_{\cA'}
$$
where $\beta$ is the natural map as in \eqref{subsec_duals}.

\begin{proposition} With the notation of
\eqref{subsec_jump-again}, suppose moreover that :
\begin{enumerate}
\alphaenu
\item
both $X$ and $X'$ are locally noetherian normal schemes
\item
$g$ maps the points of codimension one of $X'$ to
points of codimension $\leq 1$ of $X$.
\end{enumerate}
Then $\gamma$ is an isomorphism.
\end{proposition}
\begin{proof} Set $\cB:=\Img\,(\cA\to j_*j^*\cA)$ and
$\cB':=g^*\cB$. It is easily seen that the image of
$\cB'$ in $j'_*j'{}^*\cB'$ agrees with the image of
$\cA'$; in view of remark \ref{rem_taus-and-Us}(ii)
we deduce that the natural maps $\cA\to\cB$ and
$\cA'\to\cB'$ induce isomorphisms
$$
\cD^{-1}_\cA\isom\cD^{-1}_\cB
\qquad
\cD^{-1}_{\cA'}\isom\cD^{-1}_{\cB'}.
$$
We may then replace $\cA$ by $\cB$, and assume that
the natural map $\cA\to j_*j^*\cA$ is a monomorphism.
By lemmata \ref{lem_dual-is-rflx} and
\ref{lem_shortened}(ii.b), both the source and target
of $\gamma$ are reflexive $\cO_{\!X'}$-modules. Thus,
denote by $V\subset X'$ the subset of all points such
that $\gamma_x$ is an isomorphism; then $V$ is an open
subset of $X'$ (\cite[Th.4.10]{Mat}) and the inclusion
map $V\to X'$ is quasi-compact, since $X'$ is locally
noetherian. In view of proposition
\ref{prop_extend-rflx}(ii), it suffices to check that
$V$ contains all the points of codimension one of $X'$.
Thus, let $x'$ be such a point, and set $x:=g(x')$;
by assumption (b), the point  $x$ has codimension $\leq 1$
in $X$, so $\cO_{\!X,x}$ is either a field or a discrete
valuation ring and $(j_*j^*\cA)_x$ is a torsion-free
$\cO_{\!X,x}$-module of finite type. Then the same holds
for $\cA_x$, so the latter is a free $\cO_{\!X,x}$-module
of finite rank. By lemma \ref{lem_jumping}, it follows
that $\eqref{eq_diff-base-change}_{x'}$ is an isomorphism.
It is also clear that $\beta_{x'}$ is an isomorphism, and
since
$\gamma_{x'}=\eqref{eq_diff-base-change}_{x'}\circ\beta_{x'}^{-1}$,
the contention follows.
\end{proof}

\begin{remark}\label{rem_Hulot-demissionne}
(i)\ \
Let us now consider a domain $A$, and a finite, injective, and
generically \'etale ring homomorphism $f:A\to B$. As in
\eqref{sec_diagonal-idempotents}, we let $X:=\Spec\,A$, and denote
by $\cB$ the quasi-coherent $\cO_X$-algebra arising from $B$. Let
also $K$ be the field of fractions of $A$; set $B_K:=K\otimes_AB$,
and let $\bar B$ be the image of $B$ in $B_K$. We may then regard
the diagonal idempotent $e_f$ of the map $f$ as an element in
$B_K\otimes_KB_K$, and the associated trace form as a $K$-bilinear
map $B_K\otimes_KB_K\to K$, inducing an isomorphism
$\omega_f:B_K\isom B_K^\vee$. With this notation, the inverse
different $\cD^{-1}_\cB$ corresponds to
$\omega^{-1}_f(\bar B{}^\vee)\subset B_K$, and the subset
$$
\cD_\cB:=\{b\in B~|~b\cD^{-1}_\cB\subset\bar B\}
$$
is an ideal of $B$ called the {\em different of $f$}.
(If $A$ and $B$ are noetherian normal domains, then
$\cD^{-1}_\cB$ is the inverse of $\cD_\cB$ in the group
of reflexive fractional ideals of $B$ : see example
\ref{ex_integral-doms}.)

(ii)\ \
Suppose moreover that $B$ is a {\em projective $A$-module}. Then
we have a natural isomorphism $B\otimes_AB^\vee\isom\Hom_A(B,B)$,
and we denote by $\zeta_f$ the unique element of $B\otimes_AB^\vee$
corresponding to $\one_B:B\to B$ under this isomorphism. According
to \cite[Rem.4.1.17]{Ga-Ra}, we have the identity :
$$
(B_K\otimes_K\omega_f)(e_f)=\zeta_f
$$
and it follows easily that :
$$
(B\otimes_A\cD_\cB)\cdot e_f\subset\Img(B\otimes_AB\to B_K\otimes_KB_K).
$$

(iii)\ \
If $B$ is not a projective $A$-module, the inclusion of (ii)
can fail; however, if $A$ is a regular and noetherian domain,
theorem \ref{th_diagonal-idempotent} tells us that we have at least :
$$
(\cD_\cB\otimes_A\cD_\cB)\cdot e_f\subset\Img(B\otimes_AB\to B_K\otimes_KB_K).
$$
\end{remark}

\subsection{Big Cohen-Macaulay algebras}
We shall use the following auxiliary construction. For every
abelian group $(G,0,+)$, set
$$
G^\diamond:=G^\N/G^{(\N)}.
$$
Hence, the elements of $G^\diamond$ are represented by
sequences $(g_n~|~n\in\N)$ of elements of $G$, and such
a sequence represents the zero element of the quotient
group $G^\diamond$ $\Leftrightarrow$ the set
$\{n\in\N~|~g_n\neq 0\}$ is finite. Especially,
$G=0\Leftrightarrow G^\diamond=0$. Moreover, for every
group homomorphism $\phi:G\to H$, we have
$\phi^\N(G^{(\N)})\subset H^{(\N)}$, hence $\phi$ induces
a group homomorphism
$$
\phi^\diamond:G^\diamond\to H^\diamond
$$
and the rules : $G\mapsto G^\diamond$ and
$\phi:\mapsto\phi^\diamond$ define an endofunctor of the
category of abelian groups. Notice as well that
$(\phi+\psi)^\diamond=\phi^\diamond+\psi^\diamond$ for every
pair of group homomorphisms $\phi,\psi:G\to H$. Lastly,
we have an injective natural transformation that maps $G$
diagonally into $G^\diamond$
$$
\Delta_G:G\to G^\diamond
\quad :\quad
g\mapsto(g,g,\dots).
\qquad
$$

$\bullet$\ \
If $(R,+,\cdot,0,1)$ is a ring, then $R^{(\N)}$ is an ideal of
the $R$-algebra $R^\N$ (the product in the category of rings of
copies of $R$ indexed by $\N$), hence $R^\diamond$ is naturally
an $R$-algebra via the map $\Delta_R$, and for every ring
homomorphism $\phi:R\to S$, the map
$\phi^\diamond:R^\diamond\to S^\diamond$ is a ring homomorphism
such that
$$
\Delta_S\circ\phi=\phi^\diamond\circ\Delta_R.
$$

$\bullet$\ \
Moreover, if $M$ is any $R$-module, then $M^\N$ is obviously
an $R^\N$-module, and since $R^{(\N)}\cdot M^\N=M^{(\N)}$, it
follows that $M^\diamond$ is naturally an $R^\diamond$-module,
and for every $R$-linear map $f:M\to N$, the map
$f^\diamond:M^\diamond\to N^\diamond$ is $R^\diamond$-linear.

\begin{lemma}\label{lem_about-diamond}
{\em(i)}\ \
The functor $(-)^\diamond:R\Mod\to R^\diamond\Mod$ is exact
and faithful.

{\em(ii)}\ \
We have $(IM)^\diamond=IM^\diamond$ for every ideal $I\subset R$
of finite type, and every $R$-module $M$.
\end{lemma}
\begin{proof}(i): Since we have $M=0\Leftrightarrow M^\diamond=0$,
it is clear that $(-)^\diamond$ is faithful. Next, let $f:M\to N$
be an $R$-linear map, and $g_\bullet:=(g_n~|~n\in\N)\in M^\N$;
denote by $\bar g_\bullet\in M^\diamond$ the class of $g_\bullet$,
and suppose that $f^\diamond(\bar g_\bullet)=0$. This means that
there exists $k\in\N$ such that $f(g_n)=0$ for every
$n\geq k$; then clearly $\bar g_\bullet$ is in the image
of the natural map $(\Ker\,f)^\diamond\to M^\diamond$.
In particular, if $\phi$ is injective, the same holds
for $\phi^\diamond$, and it is also clear that if $\phi$
is surjective, the same holds for $\phi^\diamond$.
It follows that the natural surjection $N\to\Coker\,f$
induces a surjection $N^\diamond\to(\Coker\,f)^\diamond$,
whose kernel is $(\Img\,f)^\diamond$, and it is easily
seen that $(\Img\,f)^\diamond=\Img(f^\diamond)$; hence
$(\Coker\,f)^\diamond$ represents $\Coker(f^\diamond)$,
whence the exactness of $(-)^\diamond$.

(ii): Let $a_1,\dots,a_k$ a finite system of generators of
$I$; there follows an $R$-linear surjection
$f:M^{\oplus k}\to IM$ whose composition with the inclusion
map $j:IM\to M$ is given by the rule :
$(m_1,\dots,m_k)\mapsto a_1m_1+\cdots+a_km_k$. Hence
$j^\diamond\circ f^\diamond=(j\circ f)^\diamond:
(M^\diamond)^{\oplus k}\to M^\diamond$ is given by the same rule,
and by the foregoing $f^\diamond:(M^\diamond)^{\oplus k}\to(IM)^\diamond$
is a surjection, whence the claim.
\end{proof}

\sset\subsubsection{}\label{subsec_local-diamond}
Let now $(V,\fm)$ be a basic setup (in the sense of
\cite[\S2.1.1]{Ga-Ra}), such that $\fm$ is the filtered
union of a countable system of principal subideals.
We consider the subset
$$
\sS\subset V^\N
$$
of all $(a_n~|~n\in\N)$ such that $Va_n\subset Va_{n+1}$
for every $n\in\N$, and $\fm\subset\bigcup_{n\in\N}Va_n$.

\begin{remark}\label{rem_about-S}
(i)\ \
Since $\fm=\fm^2$, it is easily seen that $\sS$ is a
multiplicative system.

(ii)\ \
Also, for every $V$-module $M$ we have
$$
\Ker(j_M:M^\diamond\to\sS^{-1}M^\diamond)=\{x\in M^\diamond~|~\fm x=0\}
\qquad\text{and}\qquad
\fm\cdot\Coker(j_M)=0.
$$
Especially, $j_M^a:(M^\diamond)^a\to(\sS^{-1}M^\diamond)^a$
is an isomorphism of $(V,\fm)^a$-modules.

(iii)\ \
Notice that $\sS$ contains the multiplicative subset
$\sT:=(e_{k,\bullet}~|~k\in\N)$, where
$e_{k,\bullet}:=(e_{k,n}~|~n\in\N)$, with $e_{k,n}:=0$ for every
$n=0,\dots,k-1$ and $e_{k,n}:=1$ for every $n\geq k$. Notice
moreover that $e_{k,\bullet}$ is invertible in $V^\diamond$ for
every $k\in\N$, and for every $V$-module $M$ the projection
$M^\N\to M^\diamond$ factors through the localization
$M^\N\to\sT^{-1}M^\N$ and an isomorphism
$$
\sT^{-1}M^\N\isom M^\diamond.
$$
Hence, we also get an isomorphism of $V^\diamond$-modules
$$
\sS^{-1}M^\diamond\isom\sS^{-1}M^\N
\qquad
\text{for every $V$-module $M$}.
$$
\end{remark}

\begin{definition}\label{def_local-diamond}
In the situation of \eqref{subsec_local-diamond}, let also
$A\to V$ be a ring homomorphism, and $M$ any $V$-module; we
say that $M^a$ is {\em flat over $A$} (resp.
{\em faithfully flat over $A$}) if the functor
$$
A\Mod\to(V,\fm)^a\Mod
\qquad
N\mapsto M^a\otimes_AN:=(M\otimes_AN)^a
$$
is exact (resp. is exact and faithful).
\end{definition}

\begin{lemma}\label{lem_local-diamond}
With the notation of definition {\em\ref{def_local-diamond}},
the following holds :
\begin{enumerate}
\item
$M^a=0$ if and only if\/ $\sS^{-1}M^\diamond=0$.
\item
The functor $V\Mod\to V^\diamond\Mod$ : $M\mapsto\sS^{-1}M^\diamond$
factors via an exact faithful functor
$$
\sS^{-1}(-)^\diamond:V^a\Mod\to V^\diamond\Mod
\qquad
M^a\mapsto\sS^{-1}M^\diamond.
$$
\item
The unit of adjunction $\sS^{-1}M^\diamond\to(\sS^{-1}M^\diamond)^a_*$
is an isomorphism of $V$-modules.
\item
If $A$ is a coherent ring, then $M^a$ is flat over $A$ if and
only if\/ $\sS^{-1}M^\diamond$ is a flat $A$-module.
\item
If $A$ is noetherian, then $M^a$ is faithfully flat over $A$
if and only if\/ $\sS^{-1}M^\diamond$ is a faithfully flat $A$-module.
\end{enumerate}
\end{lemma}
\begin{proof} The proof of (i) shall be left to the reader.

(ii) is an immediate consequence of (i) and lemma
\ref{lem_about-diamond}(i).

(iii): Let $x\in\sS^{-1}M^\diamond$ be an element such that
$\fm\cdot x=0$, and write $x=s_\bullet^{-1}m_\bullet$ for some
$s_\bullet:=(s_n~|~n\in\N)\in\sS$ and
$m_\bullet:=(m_n~|~n\in\N)\in M^\diamond$. Let us also write
$\fm=\bigcup_{i\in\N}Va_n$ for a sequence $(a_i~|~i\in\N)$
of elements of $V$ such that $Va_i\subset\fm\cdot a_{i+1}$
for every $i\in\N$. Then, for every $i\in\N$ there exist
$\rho(i)\in\N$, $c_i\in\fm$, and
$t_{i,\bullet}:=(t_{i,n}~|~n\in\N)\in\sS$ such that
$$
a_i=c_ia_{i+1}
\qquad\text{and}\qquad
a_im_nt_{i,n}=0
\qquad
\text{for every $n\geq\rho(i)$}
$$
and we may assume that $\rho(i+1)>\rho(i)$ for every $i\in\N$.
Moreover, since $\fm\subset\bigcup_{n\in\N}Vt_{i,n}$, we may also
assume that $a_i\in Vt_{i,n}$ for every $n\geq\rho(i)$, in which
case we get :
$$
a_i^2m_n=0
\qquad
\text{for every $i\in\N$ and $n\geq\rho(i)$}.
$$
Hence, let $b_\bullet:=(b_n~|~n\in\N)\in V^\N$ be the element such
that $b_n:=0$ for every $n<\rho(0)$, and $b_n:=a^2_i$ for every
$i\in\N$ and every $n=\rho(i),\dots,\rho(i+1)-1$. It follows easily
that $b_\bullet\in\sS$ and $b_\bullet\cdot m_\bullet=0$ in $M^\diamond$,
so finally $x=0$. This shows that the unit of adjunction
$\eta_{\sS^{-1}M^\diamond}:\sS^{-1}M^\diamond\to(\sS^{-1}M^\diamond)^a_*$ is
injective. Next, let $y\in(\sS^{-1}M^\diamond)^a_*$. By remark
\ref{rem_make-comp-pair-alm-struct}(ii), the $V$-module
$\tilde\fm$ is the inductive limit of a system $(L_n~|~n\in\N)$
with $L:=V$ for every $n\in\N$, and with $V$-linear transition
maps $\phi_n:L_n\to L_{n+1}$ such that $\phi_n(1)=c_n$ for every
$n\in\N$; hence $y$ is represented by a system
$(y^{(n)}~|~n\in\N)$ of elements of $\sS^{-1}M^\diamond$ such that
\set\begin{equation}\label{eq_relation}
c_ny^{(n+1)}=y^{(n)}
\qquad
\text{for every $n\in\N$}.
\end{equation}

\begin{claim}\label{cl_simplify}
Let $\Sigma\subset\sS$ be any countable subset. Then
$\sS\cap\bigcap_{s_\bullet\in\Sigma}V^\diamond s_\bullet\neq\emptyset$.
\end{claim}
\begin{pfclaim} Say that $\Sigma=\{s_{n,\bullet}~|~n\in\N\}$.
Let the sequence $(a_n~|~n\in\N)$ of generators of $\fm$ be as
in the foregoing; then, for every $n,k\in\N$ there exists
$\rho(n,k)\in\N$ such that
$$
a_k\in Vs_{n,j}
\qquad
\text{for every $j\geq\rho(n,k)$}.
$$
Hence, for every $k\in\N$, we can find $\rho'(k)\in\N$ such that
$$
a_k\in Vs_{n,j}
\qquad
\text{for every $n=0,\dots,k$ and every $j\geq\rho'(k)$}
$$
and moreover we may assume that $\rho'(k+1)>\rho'(k)$ for every
$k\in\N$. Then, let $x_\bullet:=(x_n~|~n\in\N)\in V^\N$ be the
sequence such that $x_n:=0$ for every $n=0,\dots,\rho'(0)-1$,
and $x_n:=a_k$ whenever $\rho'(k)\leq n<\rho'(k+1)$, for
every $k\in\N$. It is easily seen that
$x_\bullet\in\sS\cap V^\diamond s_{n,\bullet}$ for every $n\in\N$,
whence the claim.
\end{pfclaim}

Now, say that $y^{(n)}=s_{n,\bullet}^{-1}m_{n,\bullet}$, with
$s_{n,\bullet}\in\sS$ and $m_{n,\bullet}\in M^\diamond$ for every
$n\in\N$; by claim \ref{cl_simplify}, we may assume that
all the elements $s_{n,\bullet}$ are equal, in which case
\eqref{eq_relation} amounts to the identities :
$$
m_{n,\bullet}=c_nm_{n+1,\bullet}
\qquad
\text{in $\sS^{-1}M^\diamond$ for every $n\in\N$}.
$$
The latter in turns means that for every $n\in\N$ there exists
$t_{n,\bullet}\in\sS$ such that
$t_{n,\bullet}m_{n,\bullet}=c_nt_{n,\bullet}m_{n+1,\bullet}$ in $M^\diamond$.
By invoking again claim \ref{cl_simplify}, we may also assume
that all the $t_{n,\bullet}$ are equal; then set
$m'_{n,\bullet}:=t_{0,\bullet}m_{n,\bullet}$ for every $n\in\N$, so that
$m'_{n,\bullet}=c_nm'_{n+1,\bullet}$ in $M^\diamond$ for every $n\in\N$.
The latter means that for every $n\in\N$ there exists $\rho(n)\in\N$
such that
$$
m_{n,j}=c_nm_{n+1,j}
\qquad
\text{for every $j\geq\rho(n)$}
$$
and we may assume that $\rho(n+1)>\rho(n)$ for every $n\in\N$.
Then, let $t_\bullet:=(t_n~|~n\in\N)\in V^\N$ (resp.
$l_\bullet:=(l_n~|~n\in\N)\in M^\N$) be the element with $t_j:=0$
(resp. $l_j:=0$) for $j=0,\dots,\rho(0)-1$ and $t_j:=a_n$ (resp.
$l_j:=m_{n,j}$) whenever $\rho(n)\leq j<\rho(n+1)$, for every $n\in\N$.

\begin{claim}\label{cl_final-detail}
$t_\bullet m_{n,\bullet}=a_nl_\bullet$ in $M^\diamond$ for every $n\in\N$.
\end{claim}
\begin{pfclaim} It suffices to show that for every $k\geq n$ and
every $j=\rho(k),\dots,\rho(k+1)-1$ we have $a_km_{n,j}=a_nm_{k,j}$.
The latter follows by a simple inspection.
\end{pfclaim}

Set $z:=t_\bullet^{-1}l_\bullet$; claim \ref{cl_final-detail} implies
that $m_{n,\bullet}=a_ny$ in $\sS^{-1}M^\diamond$ for every $n\in\N$;
hence $\eta_{\sS^{-1}M^\diamond}(z)=y$, whence the assertion.

(iv): Suppose that $M^a$ is flat over $A$, and let $I\subset A$ be
an ideal of finite type; in light of lemma \ref{lem_about-diamond}(ii),
it suffices to check that the natural map
$$
f:I\otimes_A\sS^{-1}M^\diamond\to\sS^{-1}(IM)^\diamond
$$
is an isomorphism. Now, it is easily seen that for every
$A$-module $N$ of finite type, the natural map
$N\otimes_A\sS^{-1}M^\diamond\to\sS^{-1}(N\otimes_AM)^\diamond$
is surjective; on the other hand, since $A$ is coherent,
there exists $k\in\N$ and a surjective $A$-linear map
$A^{\oplus k}\to I$, whose kernel is an $A$-module $N$ of
finite type. Thus, we have the following commutative diagram
$$
\xymatrix{ 0 \ar[r] & N\otimes_A\sS^{-1}M^\diamond \ar[r] \ar[d] &
A^{\oplus k}\otimes_A\sS^{-1}M^\diamond \ar[r] \ar[d] &
I\otimes_A\sS^{-1}M^\diamond \ar[r] \ar[d]^f & 0 \\
0 \ar[r] & \sS^{-1}(N\otimes_AM)^\diamond \ar[r] &
\sS^{-1}(M^{\oplus k})^\diamond \ar[r] &
\sS^{-1}(IM)^\diamond \ar[r] & 0
}$$
whose first and third vertical arrows are surjections, and
whose central vertical arrow is an isomorphism. Moreover,
the top horizontal sequence is right exact, and since $M^a$
is flat over $A$, combining (i) and lemma
\ref{lem_about-diamond}(i) it is easily seen that the
bottom horizontal arrow is a short exact sequence. Then,
by the snake lemma, $\Ker\,f=0$, whence the contention.

Conversely, suppose that $\sS^{-1}M^\diamond$ is a flat $A$-module,
and consider any short exact sequence
$\Sigma:=(0\to N'\to N\to N''\to 0)$ of $A$-modules; we need
to show that the induced sequence $M^a\otimes_A\Sigma$ is
still exact. To this aim, notice that $\Sigma$ is the filtered
colimit of a system of short exact sequences
$(0\to N'_\lambda\to N_\lambda\to N''_\lambda\to 0~|~\lambda\in\Lambda)$
such that $N_\lambda$ and $N''_\lambda$ are finitely presented
$A$-modules, for every $\lambda\in\Lambda$. Thus, we may assume
that $N$ and $N''$ is finitely presented, in which case the same
holds for $N'$, since $A$ is coherent, and we consider the
diagram :
$$
\xymatrix@C-10pt{ \Sigma\otimes_A\sS^{-1}M^\diamond \ar[d] &
0 \ar[r] & N'\otimes_A\sS^{-1}M^\diamond \ar[r] \ar[d] &
N\otimes_A\sS^{-1}M^\diamond \ar[r] \ar[d] &
N''\otimes_A\sS^{-1}M^\diamond \ar[r] \ar[d] & 0 \\
\sS^{-1}((\Sigma\otimes_AM)^\diamond) &
0 \ar[r] & \sS^{-1}(N'\otimes_AM)^\diamond \ar[r] &
\sS^{-1}(N\otimes_AM)^\diamond \ar[r] &
\sS^{-1}(N''\otimes_AM)^\diamond \ar[r] & 0.
}$$
Since $N$, $N'$ and $N''$ are finitely presented, lemma
\ref{lem_about-diamond}(i) easily implies that all the
vertical arrows are isomorphisms. By assumption, the top
horizontal row is short exact, hence the same holds for
the bottom one, and therefore the sequence of $(V,\fm)^a$-modules
$((\Sigma\otimes_AM)^\diamond)^a$ is short exact; then we
conclude by invoking again lemma \ref{lem_about-diamond}(i).

(v): Suppose that $M^a$ is faithfully flat over $A$; in view
of (iv), it suffices to show that for every $A$-module $N\neq 0$
of finite type, we have $\sS^{-1}M^\diamond\otimes_AN\neq 0$. But
then $N$ is finitely presented, since $A$ is noetherian, hence
the natural map
$\sS^{-1}M^\diamond\otimes_AN\to\sS^{-1}(M\otimes_AN)^\diamond$
is an isomorphism, and by assumption $M^a\otimes_AN\neq 0$;
then the assertion follows from (i). One argues similarly
for the converse : details left to the reader.
\end{proof}

\begin{proposition}\label{prop_diamond-is-complete}
In the situation of \eqref{subsec_local-diamond}, let
$I\subset V$ be any ideal of finite type; then both $M^\diamond$
and $\sS^{-1}M^\diamond$ are $I$-adically complete $V$-modules.
\end{proposition}
\begin{proof} Let $x_{\bullet\bullet}:=(x_{k,\bullet}~|~k\in\N)$
be a Cauchy sequence for the $I$-adic topology in $M^\diamond$.
Hence, $x_{k,\bullet}$ is the class of a sequence
$(x_{k,i}~|~i\in\N)$, for every $i\in\N$, and after replacing
$x_{\bullet\bullet}$ by a subsequence, we may assume that
$x_{k,\bullet}-x_{k+1,\bullet}\in I^kM^\diamond$ for every $k\in\N$.
In light of lemma \ref{lem_about-diamond}(ii), the latter
means that for every $k\in\N$ there exists $n(k)\in\N$ such
that $x_{k,i}-x_{k+1,i}\in I^kM$ for every $i\geq n(k)$, and
clearly we may assume that $n(k+1)>n(k)$ for every $k\in\N$.
Now, let $y_\bullet:=(y_i~|~i\in\N)$ be the following element
of $M^\diamond$. For every $i\leq n(0)$ we let $y_i:=x_{0,i}$;
then, for every $k\in\N$ and
$n(k)<i\leq n(k+1)$ we let $y_i:=x_{k+1,i}$. For every
$k\geq j$ and $n(k)<i\leq n(k+1)$, we get :
$$
y_i-x_{j,i}=(x_{k+1,i}-x_{k,i})+\cdots+(x_{j+1,i}-x_{j,i})\in
I^kM+\cdots+I^jM=I^jM.
$$
Hence, $y_\bullet-x_{j,\bullet}\in I^jM^\diamond$ for every
$j\in\N$, so $y_\bullet$ is a limit of $x_{\bullet\bullet}$.

Next, let $x_{\bullet\bullet}:=(x_{k,\bullet}~|~k\in\N)$ be a Cauchy
sequence for the $I$-adic topology in $\sS^{-1}M^\diamond$; again,
we may assume that
$y_{k,\bullet}:=x_{k+1,\bullet}-x_{k,\bullet}\in I^k\sS^{-1}M^\diamond$
for every $k\in\N$, and then there exist $s_{k,\bullet}\in\sS$
and $z_{k,\bullet}\in I^kM^\diamond$ with
$y_{k,\bullet}=s_{k,\bullet}^{-1}\cdot z_{k,\bullet}$ for every
$k\in\N$. Also, we may assume :
\set\begin{equation}\label{eq_divisibility}
s_{k+1,\bullet}\in V^\N s_{k,\bullet}
\qquad
\text{for every $k\in\N$}.
\end{equation}
Now, by assumption there exists a sequence $(a_n~|~n\in\N)$
of elements of $V$ such that $Va_n\subset Va_{n+1}$ for every
$n\in\N$, and $\bigcup_{n\in\N}Va_n=\fm$. For every
$k\in\N$ we define inductively a sequence of integers
$(n(k,i)~|~i\in\N)$ as follows. For $i=0$, let $n(k,0)$
be the smallest of the integers $n\in\N$ such that
$a_0\in Vs_{k,n}$; then, for every $i>0$ let $n(k,i)\in\N$
be the smallest of the integers $n>n(k,i-1)$ such that
$a_i\in Vs_{k,n}$. Notice that \eqref{eq_divisibility} implies
$$
n(k+1,i)\geq n(k,i)
\qquad
\text{for every $k,i\in\N$}.
$$
Lastly, let $t_\bullet$ be the sequence such that $t_n:=0$ for
every $n<n(0,0)$, and $t_n:=a_i$ for every $n,i\in\N$ with
$n(i,i)\leq n<n(i+1,i+1)$. We have $t_\bullet\in\sS$, and
$t_i\in Vs_{k,i}$ for every $k\in\N$ and every $i\geq n(k,k)$;
thus, the class of $s_{k,\bullet}$ in $V^\diamond$ divides the class
of $t_\bullet$, and we may replace $s_{k,\bullet}$ by $t_\bullet$ for
every $k\in\N$. By the foregoing case, the series
$\sum_{k\in\N}z_{k,\bullet}$ converges in $M^\diamond$, and let
$z_\bullet$ be a limit; then $t_\bullet^{-1}z_\bullet$ is a limit
for $x_{\bullet\bullet}$ in $\sS^{-1}M^\diamond$.
\end{proof}

\begin{lemma}\label{lem_almost-perfectoid}
In the situation of \eqref{subsec_almost-perfectoids}, let us
endow $\sS^{-1}A^\diamond_0$ with its $b$-adic topology, and let
$A_0^\wedge$ be the $b$-adic completion of $A_0$; then the
basic setup $(A^\wedge_0,\fm A^\wedge_0)$ is almost perfectoid if
and only if the maximal separated quotient
$(\sS^{-1}A^\diamond_0)^\sep$ of\/ $\sS^{-1}A^\diamond_0$ is perfectoid.
\end{lemma}
\begin{proof} Clearly
$\sS^{-1}\bar\Phi{}_{A_0}^\diamond:\sS^{-1}(A_0/bA_0)^\diamond\to
\sS^{-1}(A_0/b^pA_0)^\diamond$ is the ring endomorphism induced
by Frobenius endomorphism of $\sS^{-1}(A_0/b^pA_0)^\diamond$. On
the other hand, by lemma \ref{lem_local-diamond}(ii),
$\bar\Phi{}^a_{A_0}$ is an isomorphism if and only if the same
holds for $\sS^{-1}\bar\Phi{}_{A_0}^\diamond$. Notice also that
$C:=(\sS^{-1}A_0^\diamond)^\sep$ is $b$-adically complete and
separated, by virtue of propositions \ref{prop_diamond-is-complete}
and \ref{prop_replaces-Mat-Th.8.1}(v).
Then let $\bar\Phi_C:C/bC\to C/b^pC$ be the ring homomorphism induced
by the Frobenius endomorphism of $C/b^pC$; in view of theorem
\ref{th_regular-seq-criterion}, we are reduced to showing :

\begin{claim} $\bar\Phi_C$ is an isomorphism if and only if the
same holds for $\sS^{-1}\bar\Phi{}_{A_0}^\diamond$.
\end{claim}
\begin{pfclaim}[] By lemma \ref{lem_about-diamond}(ii),
the projection $A_0\to A_0/bA_0$ induces a ring isomorphism
$\sS^{-1}A^\diamond_0/(b)\isom\sS^{-1}(A_0/bA_0)^\diamond$, and
likewise $\sS^{-1}A^\diamond_0/(b^p)\isom\sS^{-1}(A_0/b^pA_0)^\diamond$.
Therefore, the projections $\sS^{-1}A_0^\diamond\to C/bC$ and
$\sS^{-1}A_0^\diamond\to C/b^pC$ factor through ring isomorphisms
$$
\sS^{-1}(A_0/bA_0)^\diamond\isom C/bC
\qquad
\sS^{-1}(A_0/b^pA_0)^\diamond\isom C/b^pC
$$
whence the claim.
\end{pfclaim}
\end{proof}

\sset\subsubsection{}\label{subsec_long-situation}
In the situation of \eqref{subsec_almost-perfectoids}, let
$g_\bullet:=(g_n~|~n\in\N)\in\bE(A_0)$, and set $g:=g_0\in A_0$,
so that $g^\gamma$ is well defined in $A_0$ for every
$\gamma\in\N[1/p]$. We suppose that
$\fm=\bigcup_{\gamma\in\N[1/p]\setminus\{0\}}A_0g^\gamma$, and
that the basic setup $(A^\wedge_0,\fm A_0^\wedge)$ is almost
perfectoid (see definition \ref{def_almost-perfectoid}).

\begin{proposition}\label{prop_massive-attack}
In the situation of \eqref{subsec_long-situation}, suppose
moreover that $A_0$ is complete and separated, and let $A_0^\nu$
be the integral closure of the image of $A_0$ in $A_0[1/g]$. Then
the induced morphism $A^a_0\to(A_0^\nu)^a$ is an isomorphism
of $(A_0,\fm)^a$-algebras.
\end{proposition}
\begin{proof} Let us check first that the kernel of the
localization $A_0\to A_0[1/g]$ is annihilated by $\fm$.
Hence, let $x\in A_0$ such that $g^nx=0$ for some $n\in\N$,
and denote by $\bar x_\bullet\in B:=(\sS^{-1}A_0^\diamond)^\sep$
the image of the constant sequence $(x,x,x,\dots)\in A_0^\N$
(notation of lemma \ref{lem_almost-perfectoid}). Since
$B$ is reduced (corollary \ref{cor_perf-are-reduced}), we
have $g^\gamma\cdot\bar x=0$ in $B$ for every
$\gamma\in\N[1/p]\setminus\{0\}$. The latter means that for
all $i\in\N$ and $\gamma\in\N[1/p]\setminus\{0\}$ there
exists $s_\bullet:=(s_n~|~n\in\N)\in\sS$ such that
$s_\bullet\cdot g^\gamma\cdot\bar x\in b^iA_0^\diamond$. In turn,
this means that for every such $i$ and $\gamma$ there exists
$s_\bullet$ and $n\in\N$ such that $s_mg^\gamma x\in b^iA_0$
for every $m\geq n$; equivalently, for every $i\in\N$ and
$\gamma\in\N[1/p]\setminus\{0\}$ we have $g^\gamma x\in b^iA_0$,
and then $g^\gamma x=0$ in $A_0$ for every such $\gamma$,
as required.

Next, we show that the cokernel of the map $A_0\to A_0^\nu$
is annihilated by $\fm$. Indeed, let $a/g^n\in A_0^\nu$
for some $a\in A_0$ and $n\in\N$, and denote by
$\bar a_\bullet\in B$ the image of the constant sequence
$(a,a,a,\dots)\in A_0^\N$. Hence $\bar a_\bullet/g^n\in B[g^{-1}]$
lies in the integral closure $C$ of the image $\bar B$ of $B$
in $B[1/g]$. Since $B$ is perfectoid, the induced map
$B^a\to C^a$ is an isomorphism of $(V,\fm)^a$-algebras
(theorem \ref{th_perf-Abhyankar}(iv)), hence
$g^\gamma\cdot(\bar a_\bullet/g^n)\in\bar B$ for every
$\gamma\in\N[1/p]\setminus\{0\}$; then there exists
$y_\bullet:=(y_n~|~n\in\N)\in B$ such that
$$
g^{\gamma+\gamma'}\cdot\bar a_\bullet=g^{n+\gamma'}\cdot y_\bullet
\qquad
\text{for every $\gamma,\gamma'\in\N[1/p]$}.
$$
Hence, for every $i\in\N$ there exists $s_\bullet\in\sS$ such that
$s_\bullet(g^{\gamma+\gamma'}\cdot\bar a_\bullet-g^{n+\gamma'}\cdot y_\bullet)
\in b^iA_0^\diamond$. The latter implies that for every $i\in\N$
and every $\gamma\in\N[1/p]\setminus\{0\}$ the image $z_{\gamma,i}$
of $g^\gamma a$ in $A_0/b^iA_0$ lies in $J_i:=g^n(A_0/b^iA_0)$.
For every $i\in\N$, let $K_i:=\Ker\,g^n\one_{A_0/b^iA_0}$; we get
an exact sequence
$$
K:=\lim_{i\in\N}K_i\to A_0=\lim_{i\in\N}A_0/b^iA_0\xrightarrow{\alpha}
J:=\lim_{i\in\N}J_i\to\lim_{n\in\N}{}^{\!1}K_i
$$
and notice that $K=\Ker\,g^n\one_{A_0}$ lies in the kernel of
the localization $A_0\to A_0[1/g]$, hence $K^a=0$, by the
foregoing. Moreover, the system of inclusions
$(J_i\to A_0/b^iA_0~|~i\in\N)$ induces a map $\beta:J\to A_0$
such that $\beta\circ\alpha=g^n\one_{A_0}$.

\begin{claim}\label{cl_take-5}
The inverse system $(K_i~|~i\in\N)$ is almost essentially zero.
\end{claim}
\begin{pfclaim} Since $A_0$ is almost perfectoid, the Frobenius
endomorphism of $A_0/b^pA_0$ induces a map $A_0/bA_0\to A_0/b^pA_0$
whose kernel is almost zero. By a simple induction on $j$, it
follows that for every $j\in\N$ and every $x\in A_0$ such that
$x^p\in b^{pj}A_0$, we have $g^\gamma x\in b^jA_0$ for every
$\gamma\in\N[1/p]\setminus\{0\}$. Now, fix $i,k\in\N$ and let
$x\in A_0$ whose class in $A_0/b^{p^ki}A_0$ lies in $K_{p^ki}$,
{\em i.e.} $g^nx\in b^{p^ki}A_0$. Then $(g^{n/p^k}x)^{p^k}\in b^{p^ki}A_0$,
and an easy induction on $j$ yields
$$
(g^{(j+1)n/p^k}x)^{p^{k-j}i}\in b^{p^{k-j}i}A_0
\qquad
\text{for every $j=0,\dots,k$}.
$$
For $j=k$, it follows that $g^{(k+1)n/p^k}x\in b^iA_0$. This shows
that for every $i\in\N$ and every $\beta\in\N[1/p]\setminus\{0\}$
there exists $k\in\N$ such that the image of $K_{p^ki}$ in $K_i$
is annhilated by $g^\beta$, whence the contention.
\end{pfclaim}

By proposition \ref{prop_last-ML-brick}(iii) and claim
\ref{cl_take-5} it follows that $\lim_{i\in\N}^1K^a_i=0$,
hence $\alpha$ is almost surjective. By construction,
$z_{\gamma,\bullet}:=(z_{\gamma,i}~|~i\in\N)\in\lim_{i\in\N}J_i$;
therefore we find for every $\gamma\in\N[1/p]\setminus\{0\}$
an element $a'\in A$ such that
$\alpha(a')=g^\gamma z_{\gamma,\bullet}$, whence
$g^na'=\beta\circ(a')=g^\gamma\beta(z_{\gamma,\bullet})=g^{2\gamma}a$.
The image of $a'$ in $A_0^\nu$ therefore equals
$g^{2\gamma}(a/g^n)$, and the proof is complete.
\end{proof}

\begin{proposition}\label{prop_general-zimbabwe}
In the situation of \eqref{subsec_long-situation}, let $h\in A_0$,
and set $B:=A_0[h^{1/p^\infty}]=A_0[X^{1/p^\infty}]/(X-h)$. Let moreover
$C$ be the $p$-integral closure of $B$ in $B[1/b]$, and $C^\wedge$
the $b$-adic completion of $C$. Then $C^a$ is a faithfully flat
$(A_0,\fm)^a$-algebra, and the basic setup $(C^\wedge,\fm C^\wedge)$
is almost perfectoid for the $b$-adic topology of $C^\wedge$.
\end{proposition}
\begin{proof} Let us regard $A^\N$ (resp. $B^\N$) as an $A$-algebra
(resp. a $B$-algebra) via the diagonal map $\Delta_A:A\to A^\N$ (resp.
$\Delta_B:B\to B^\N$), and denote by $D$ the $p$-integral closure
of $B^\N$ in $B^\N[1/b]$. By remark \ref{rem_charact-p-int-clos}(ii),
the map $\Delta_B$ restricts to a ring homomorphism $j_C:C\to D$.
By the same token, for every $n\in\N$, the projection $B^\N\to B$
on the $n$-th factor of $B^\N$ restricts to a map $\pi_{C,n}:D\to C$
such that $\pi_{C,n}\circ j_C=\one_C$; let also $\pi_{A,n}:A^\N\to A$
be the corresponding projection on the $n$-th factor, and likewise
define $\pi_{M,n}:M^\N\to M$ and the diagonal map $\Delta_M:M\to M^\N$,
for every $A$-module $M$. For every such $M$ and every $n\in\N$ we
then get maps of $A_0$-modules :
$$
T_0:=\Tor_1^{A_0}(M,C)\xrightarrow{\alpha}T:=\Tor_1^{A_0^\N}(M^\N,D)
\xrightarrow{\beta_n}T_0
$$
where $\alpha$ is induced by $(\Delta_A,\Delta_M,j_C)$, and $\beta_n$
is induced by $(\pi_{A,n},\pi_{M,n},\pi_{C,n})$. On the other hand,
$\sS^{-1}A_0^\N$ is a formal perfectoid ring for its $b$-adic
topology, by lemma \ref{lem_almost-perfectoid} and remark
\ref{rem_about-S}(iii), hence $\sS^{-1}D$ is a faithfully flat
$\sS^{-1}A^\N$-algebra, by theorem \ref{th_Zimbabwe}(i) and
lemma \ref{lem_crit-p-integr-closed}(iv). Hence, $\sS^{-1}T=0$.
Now, let $x\in T_0$; then there exists
$s_\bullet:=(s_n~|~n\in\N)\in\sS$ with $s_\bullet\cdot\alpha(x)=0$
in $T$, whence $s_n\cdot\beta_n\circ\alpha(x)=s_nx=0$ in $T_0$
for every $n\in\N$, so $T_0^a=0$. This shows that $C^a$ is a
flat $A_0^a$-algebra. Next, consider, for every $n\in\N$, the
commutative diagram :
$$
\xymatrix@C+60pt{ M \ar[r]^-{\Delta_M} \ar[d]_\gamma &
M^\N \ar[d]^\delta \\
C\otimes_{A_0}M \ar[r]^-{j_C\otimes_{A_0}\Delta_M} &
D\otimes_{A_0^\N}M^\N
}$$
and let $y\in\Ker\,\gamma$; then $\Delta_M(y)\in\Ker\,\delta$,
and since $\sS^{-1}D$ is a faithfully flat $\sS^{-1}A_0^\N$-algebra,
it follows that there exists $s_\bullet\in\sS$ with
$s_\bullet\cdot\Delta_M(y)=0$, {\em i.e.} $s_ny=0$ for every
$n\in\N$, so that $(\Ker\,\gamma)^a=0$, which proves that
$C^a$ is a faithfully flat $A^a_0$-algebra.

Next, since $b$ is a regular element of $C$, it is also regular
in $C^\wedge$; also, we already know that the Frobenius endomorphism
of $C^\wedge$ induces an injection $\bar\Phi_C:C/bC\to C/b^pC$
(lemma \ref{lem_crit-p-integr-closed}(i)); on the other
hand, for any $n\in\N$, consider the commutative diagram :
$$
\xymatrix@C+60pt{
C/bC \ar[r]^-{A_0/bA_0\otimes_{A_0}j_C} \ar[d]_{\bar\Phi_C} &
D/bD \ar[r]^-{A_0/bA\otimes_{A_0}\pi_{C,n}}
\ar[d]_{\bar\Phi_D} & C/bC \ar[d]^{\bar\Phi_C} \\
C/b^pC \ar[r]^-{A_0/b^pA_0\otimes_{A_0}j_C} &
D/b^pD \ar[r]^-{A_0/b^pA\otimes_{A_0}\pi_{C,n}} & C/b^pC.
}$$
By theorem \ref{th_Zimbabwe}(ii), the map $\sS^{-1}\bar\Phi_D$ is
an isomorphism. Now, let $x\in C/b^pC$, and pick
$s_\bullet:=(s_n~|~n\in\N)\in\sS$ and $y_\bullet\in D/bD$ such that
$\sS^{-1}\bar\Phi_D(s_\bullet^{-1}y_\bullet)=(A_0/b^pA_0\otimes_{A_0}j_C)(x)$.
Set $z_n:=(A_0/bA\otimes_{A_0}\pi_{C,n})(y_\bullet)$ for every
$n\in\N$; then $\bar\Phi_C(z_n)=s_nx$ for every $n\in\N$.
This shows that $\bar\Phi_C$ is an almost isomorphism, so
$(C^\wedge,\fm C^\wedge)$ is an almost perfectoid basic setup.
\end{proof}

\begin{remark}\label{rem_alternative}
(i)\ \
In the situation of \eqref{subsec_long-situation}, suppose
moreover that there exists a perfectoid Tate ring $R$, with
a perfectoid subring of definition $R_0$, and an adic ring
homomorphism $R_0\to A_0$, such that $g_\bullet$ is the image of
an element $g'_\bullet:=(g'_n~|~n\in\N)$ of $\bE(R_0)$. Then we
get a basic setup $(R_0,\fm_R)$ with
$\fm_R:=\bigcup_{\gamma\in\N[1/p]\setminus\{0\}}g'^\gamma_0R_0$, and we can
form as in \eqref{subsec_adjoints-for-PerfTate} the topological
$R$-algebra $A_0^\natural:=R\otimes_{W(\bE_0(R_0))}W(\bE(A_0/pA_0))$,
which is perfectoid and admits a natural adic map
$A^\natural_0\to A_0$; the latter is an almost isomorphism relative
to the basic setup $(R_0,\fm_R)$, by theorem
\ref{th_about-Andre-construction}(iii.b). Let $A_0^\nu$ (resp.
$(A_0^\natural)^\nu$) be the integral closure of $A_0$ in $A_0[1/g]$
(resp. of $A_0^\natural$ in $A_0^\natural[1/g'_0]$). Then, on the one
hand, the induced morphism $((A_0^\natural)^\nu)^a\to(A_0^\nu)^a$ is
an isomorphism of $(R_0,\fm_R)^a$-algebras (\cite[Lemma 8.2.28]{Ga-Ra});
on the other hand, the same holds for the morphism
$(A_0^\natural)^a\to((A_0^\natural)^\nu)^a$, due to theorem
\ref{th_perf-Abhyankar}(iv). In this way we thus obtain a shorter
proof of proposition \ref{prop_massive-attack}, under the stated
further assumptions.

(ii)\ \
We can also give an alternative proof of the first assertion
of proposition \ref{prop_general-zimbabwe} along the same lines,
under the assumptions of (i). Indeed, first let us remark that,
since the map $A_0^\natural\to A_0$ is an almost isomorphism, there
exists an almost element $l\in(A_0^\natural)^a_*$ whose image
in $(A_0)^a_*$ agrees with $h$; by definition,
$l:\fm_R\otimes_R\fm_R\to A_0^\natural$ is an $R$-linear map,
and we set $l_{\gamma+\delta}:=l(g'^\gamma_0\otimes g'^\delta_0)$ for
every $\gamma,\delta\in\N[1/p]\setminus\{0\}$. It is easily
seen that $l_{\gamma+\delta}$ depends only on the sum $\gamma+\delta$,
hence we get a compatible system
$(l_\gamma~|~\gamma\in\N[1/p]\setminus\{0\})$ of elements of
$A_0^\natural$ such that $g'^\delta_0\cdot l_\gamma=l_{\gamma+\delta}$
for every such $\gamma$ and $\delta$. Let $D_\gamma$ (resp.
$E_\gamma$) be the $p$-integral closure of
$A_{0,\gamma}^\natural:=A_0^\natural[l_\gamma^{1/p^\infty}]$ in
$A^\natural_{0,\gamma}[1/b]$ (resp. of
$A_{0,\gamma}:=A_0[(g^\gamma h)^{1/p^\infty}]$ in $A_{0,\gamma}[1/b]$) for
every $\gamma\in\N[1/p]\setminus\{0\}$; by theorem \ref{th_Zimbabwe}(i),
$D_\gamma$ is a faithfully flat $A_0^\natural$-algebra; moreover, for
every $\gamma,\delta\in\N[1/p]\setminus\{0\}$, the natural map
$A^\natural_{0,\gamma+\delta}\to A^\natural_{0,\gamma}$ induces a map of
$A_0^\natural$-algebras $D_{\gamma+\delta}\to D_\gamma$ (remark
\ref{rem_charact-p-int-clos}(ii)). Furthermore, for every such
$\gamma$ we have a natural map
$i_\gamma:A^\natural_{0,\gamma}\to A_{0,\gamma}$, whence an induced map
$j_\gamma:D_\gamma\to E_\gamma$. Since $p\in b^pA_0$, the localization
$A_{0,\gamma}\to A_{0,\gamma}[1/b]$ induces an isomorphism
$A_{0,\gamma}[1/p]\isom A_{0,\gamma}[1/b,1/p]$, and therefore
$E_\gamma$ is also the weak normalization of $A_{0,\gamma}$ in
$A_{0,\gamma}[1/b]$ (proposition \ref{prop_many-rings-equal});
likewise we see that $D_\gamma$ is the weak normalization of
$A^\natural_{0,\gamma}$ in $A^\natural_{0,\gamma}[1/b]$. On the other
hand, clearly $i_\gamma^a$ is an isomorphism; taking into account
\ref{rem_alm-int-closure}(iii), we conclude that the same holds
for $j_\gamma^a$. Lastly, the compatible system of ring
homomorphisms $(A_0[(g^\gamma h)^{1/p^\infty}]\to A_0[h^{1/p^\infty}]~|~
\gamma\in\N[1/p]\setminus\{0\})$ yields a compatible system
of maps $(E_\gamma\to B~|~\gamma\in\N[1/p]\setminus\{0\})$, and
it is easily seen that the induced co-cone of $A_0^a$-algebras
$(E^a_\gamma\to B^a~|~\gamma\in\N[1/p]\setminus\{0\})$ is universal.
Summing up, $B^a$ is the colimit of the filtered system of
faithfully flat $A_0^a$-algebras
$(D_\gamma^a~|~\gamma\in\N[1/p]\setminus\{0\})$, so it is
a faithfully flat $A_0^a$-algebra (lemma
\ref{lem_alm-faithful-flatness}(ii)).

(iii)\ \
Likewise, by applying the functor $(-)^\natural$ in a
similar fashion, one may show that the basic setup
$(C^\wedge,\fm C^\wedge)$ is almost perfectoid : we
leave the details as an exercise for the reader.
\end{remark}


\printindex


\begin{thebibliography}{99}

\bibitem[1]{AJL}
{\scshape L.Alonso Tarr\'{\i}o, A.Jerem\'{\i}as L\'opez, J.Lipman},
Local homology and cohomology on schemes.
{\it Ann. Sci. E.N.S.} 30 (1997) pp.1--39.

\bibitem[2]{An}
{\scshape M.Andr\'e},
Homologie des alg{\`e}bres commutatives.
{\it Springer Grundl. Math. Wiss.} 206 (1974).

\bibitem[3]{AnII}
{\scshape M.Andr\'e},
Localisation de la lissit\'e formelle.
{\it Manuscr. Math.} 13 (1974) pp.297--307.

\bibitem[4]{YAn}
{\scshape Y.Andr\'e},
La conjecture du facteur direct.
{\it \it Publ. Math. IHES} 127 (2018) pp.71--93.

\bibitem[5]{YAn-2}
{\scshape Y.Andr\'e},
Le lemme d'Abhyankar perfecto\"ide.
{\it \it Publ. Math. IHES} 127 (2018) pp.1--70.

\bibitem[6]{SGA4-1}
{\scshape M.Artin et al.},
Th\'eorie des topos et cohomologie \'etale des sch\'emas (SGA4) -- tome 1.
{\it Springer Lect. Notes Math.} 269 (1972).

\bibitem[7]{SGA4-2}
{\scshape M.Artin et al.},
Th\'eorie des topos et cohomologie \'etale des sch\'emas (SGA4) -- tome 2.
{\it Springer Lect. Notes Math.} 270 (1972).

\bibitem[8]{SGA4-3}
{\scshape M.Artin et al.}, Th{\'e}orie des topos et cohomologie
{\'e}tale des sch\'emas -- tome 3.
{\it Springer Lect. Notes Math.} 305 (1973).

\bibitem[9]{BaMu}
{\scshape H.Bass, M.P.Murthy},
Grothendieck groups and Picard groups of abelian group rings.
{\it Ann. of Math.} 86 (1967) pp.16--73.

\bibitem[10]{Bei}
{\scshape A.Beilinson},
$p$-adic periods and derived deRham cohomology.
{\it J. Amer. Math. Soc.} 25 (2012) pp.715--738.

\bibitem[11]{BBD}
{\scshape A.Beilinson, J.Bernstein, P.Deligne, O.Gabber},
Faisceaux pervers (new edition).
{\it Ast\'erisque} 100 (2018).

\bibitem[12]{SGA6}
{\scshape P.Berthelot et al.},
Th\'eorie des Intersection et Th\'eor\`emes de Riemann-Roch.
{\it Springer Lect. Notes Math.} 225 (1971).

\bibitem[13]{Bha}
{\scshape B.Bhatt},
Almost direct summands.
{\it Nagoya Math. J.} 214 (2014) pp.195--204.

\bibitem[14]{Bha2}
{\scshape B.Bhatt},
On the direct summand conjecture and its derived variant.
{\it Invent. Math.} 212 (2018) pp.297--317.

\bibitem[15]{Bj}
{\scshape J.-E. Bj\"ork},
Analytic $\cD$-Modules and Applications.
{\it Kluwer Math. and Its Appl.} 247 (1993).

\bibitem[16]{Bor}
{\scshape F.Borceux},
Handbook of Categorical Algebra I. Basic Category Theory.
{\it Cambridge Univ. Press Encycl. of Math. and Its Appl.} 50 (1994).

\bibitem[17]{BorII}
{\scshape F.Borceux},
Handbook of Categorical Algebra II. Categories and structures.
{\it Cambridge Univ. Press Encycl. of Math. and Its Appl.} 51 (1994).

\bibitem[18]{BorIII}
{\scshape F.Borceux},
Handbook of Categorical Algebra III. Categories of Sheaves.
{\it Cambridge Univ. Press Encycl. of Math. and Its Appl.} 52 (1994).

\bibitem[19]{Borg}
{\scshape J.Borger},
The basic geometry of Witt vectors -- I.
{\it Algebra and Number Theory} 5 (2011) pp.231--285.

\bibitem[20]{Bo-Ray}
{\scshape S.Bosch, W.L\"utkebohmert, M.Raynaud},
Neron Models.
{\it Springer Ergebn. der Math.} 21 (1990).

\bibitem[21]{Bo-Gun}
{\scshape S.Bosch, U.G\"untzer, R.Remmert},
Non-Archimedean analysis.
{\it Springer Grundl. Math. Wiss.} 261 (1984).

\bibitem[22]{BouAC}
{\scshape N.Bourbaki},
Alg\`ebre Commutative -- Chapitres 1 \`a 9.
{\it Hermann} (1961).

\bibitem[23]{BouAC10}
{\scshape N.Bourbaki},
Alg\`ebre Commutative -- Chapitre 10.
{\it Masson} (1998).

\bibitem[24]{Bou-Ens}
{\scshape N.Bourbaki},
Th\'eorie des Ensembles.
{\it Hermann} (1963).

\bibitem[25]{BouAH}
{\scshape N.Bourbaki},
Alg\`ebre Homologique.
{\it Masson} (1980).

\bibitem[26]{Bourbaki}
{\scshape N.Bourbaki},
Alg{\`e}bre.
{\it Hermann} (1970).

\bibitem[27]{BouEVT}
{\scshape N.Bourbaki},
Espaces Vectoriels Topologiques.
{\it Hermann} (1966).

\bibitem[28]{BouTG}
{\scshape N.Bourbaki},
Topologie G{\'e}n{\'e}rale.
{\it Hermann} (1971).

\bibitem[29]{Brod}
{\scshape M.Brodmann},
Asymptotic stability of $\Ass(M/I^nM)$.
{\it Proc. Amer. Math. Soc.} 74 (1979) pp.16--18.

\bibitem[30]{Br-He}
{\scshape W.Bruns, J.Herzog},
Cohen-Macaulay rings -- revised edition.
{\it Cambridge Univ. Press} (1998).

\bibitem[31]{Bu-Ve}
{\scshape K.Buzzard, A.Verberkmoes},
Stably uniform affinoids are sheafy.
{\it J. Reine Angew. Math.} 740 (2018) pp.25--40.

\bibitem[32]{Ch}
{\scshape L.G.Chouinard II},
Krull semigroups and divisor class groups.
{\it Can. J. Math.} 33 (1981) pp.1459--1468.

\bibitem[33]{Con2}
{\scshape B.Conrad},
Deligne's notes on Nagata compactifications.
{\it J. Ramanujan Math. Soc.} 22 (2007) pp.205--257; erratum
24 (2009) pp.427--428.

\bibitem[34]{Con}
{\scshape B.Conrad},
Grothendieck Duality and Base Change.
{\it Springer Lect. Notes Math.} 1750 (2000).

\bibitem[35]{SGA3}
{\scshape M.Demazure, A.Grothendieck et al.},
Sch\'emas en Groupes I.
{\it Springer Lect. Notes Math.} 151 (1970).

\bibitem[36]{EGAI-new}
{\scshape J.Dieudonn\'e, A.Grothendieck},
\'El\'ements de G\'eom\'etrie Alg\'ebrique -- Chapitre I.
{\it Springer Grundl. Math. Wiss.} 166 (1971).

\bibitem[37]{EGAI}
{\scshape J.Dieudonn\'e, A.Grothendieck},
\'El\'ements de G\'eom\'etrie Alg\'ebrique -- Chapitre I.
{\it Publ. Math. IHES} 4 (1960).

\bibitem[38]{EGAII}
{\scshape J.Dieudonn\'e, A.Grothendieck},
\'El\'ements de G\'eom\'etrie Alg\'ebrique -- Chapitre II.
{\it Publ. Math. IHES} 8 (1961).

\bibitem[39]{EGAIII}
{\scshape J.Dieudonn\'e, A.Grothendieck},
\'El\'ements de G\'eom\'etrie Alg\'ebrique -- Chapitre III, partie 1.
{\it Publ. Math. IHES} 11 (1961).

\bibitem[40]{EGAIII-2}
{\scshape J.Dieudonn\'e, A.Grothendieck},
\'El\'ements de G\'eom\'etrie Alg\'ebrique -- Chapitre III, partie 2.
{\it Publ. Math. IHES} 17 (1963).

\bibitem[41]{EGAIV}
{\scshape J.Dieudonn\'e, A.Grothendieck},
\'El\'ements de G\'eom\'etrie Alg\'ebrique -- Chapitre IV, partie 1.
{\it Publ. Math. IHES} 20 (1964).

\bibitem[42]{EGAIV-2}
{\scshape J.Dieudonn\'e, A.Grothendieck},
\'El\'ements de G\'eom\'etrie Alg\'ebrique --
Chapitre IV, partie 2.
{\it Publ. Math. IHES} 24 (1965).

\bibitem[43]{EGAIV-3}
{\scshape J.Dieudonn\'e, A.Grothendieck},
\'El\'ements de G\'eom\'etrie Alg{\'e}brique --
Chapitre IV, partie 3.
{\it Publ. Math. IHES} 28 (1966).

\bibitem[44]{EGA4}
{\scshape J.Dieudonn\'e, A.Grothendieck},
\'El\'ements de G\'eom{\'e}trie Alg\'ebrique -- Chapitre IV, partie 4.
{\it Publ. Math. IHES} 32 (1967).

\bibitem[45]{Dut}
{\scshape P.Dutton},
Prime ideals attached to a module.
{\it Quart. J. Math. Oxford} 29 (1978) pp.403--413. 

\bibitem[46]{Ep-Sha}
{\scshape N.Epstein, J.Shapiro},
Strong Krull primes and flat modules.
{\it J. Pure Appl. Alg.} 218 (2014) pp.1712-1729.

\bibitem[47]{Fa2}
{\scshape G.Faltings},
Almost {\'e}tale extensions.
{\it Ast{\'e}risque} 279 (2002) pp.185--270.

\bibitem[48]{Fu-Ga-Ka}
{\scshape K.Fujiwara, O.Gabber, F.Kato},
On Hausdorff completions of commutative rings in rigid geometry.
{\it J. of Algebra} 332 (2011) pp.293--321.

\bibitem[49]{Fu-Ka}
{\scshape K.Fujiwara, F.Kato},
Foundations of Rigid Geometry I.
{\it EMS Monographs in Math., European Math. Society} (2018).

\bibitem[50]{Ful}
{\scshape W.Fulton},
Introduction to Toric Varieties.
{\it Princeton Univ. Press Ann. of Math. Studies} 131 (1993).

\bibitem[51]{Gab}
{\scshape O.Gabber},
Affine analog of the proper base change theorem.
{\it Isr. J. of Math.} 87 (1994) pp.325--335.

\bibitem[52]{Ga-Ra}
{\scshape O.Gabber, L.Ramero},
Almost ring theory.
{\it Springer Lect. Notes Math.} 1800 (2003).

\bibitem[53]{Ga}
{\scshape P.Gabriel},
Des cat{\'e}gories ab{\'e}liennes.
{\it Bull. Soc. Math. France} 90 (1962) pp.323--449.

\bibitem[54]{Gi}
{\scshape J.Giraud},
Cohomologie non ab{\'e}lienne.
{\it Springer Grundl. Math. Wiss.} 179 (1971).

\bibitem[55]{God}
{\scshape R.Godement},
Th{\'e}orie des faisceaux.
{\it Hermann} (1958).

\bibitem[56]{Go-Ja}
{\scshape P.G.Goerss, J.F.Jardine},
Simplicial Homotopy Theory.
{\it Birkh\"auser Modern Classics} (2009).

\bibitem[57]{Gray}
{\scshape J.W.Gray},
Formal Category Theory : Adjointness for $2$-Categories.
{\it Springer Lecture Notes Math.} 391 (1974).

\bibitem[58]{SGA1}
{\scshape A.Grothendieck et al.},
Rev{\^e}tements {\'E}tales et Groupe Fondamental.
{\it Documents Math. Soc. Math. France} 3 (2003).

\bibitem[59]{LCoh}
{\scshape A.Grothendieck},
Local cohomology.
{\it Springer Lect. Notes Math.} 41 (1967).

\bibitem[60]{Toh}
{\scshape A.Grothendieck},
Sur quelques points d'alg{\'e}bre homologique.
{\it Tohoku Math.J.} 9 (1957) pp.119-221.

\bibitem[61]{SGA2}
{\scshape A.Grothendieck},
Cohomologie locale des faisceaux coh{\'e}rents et th{\'e}or{\'e}mes de
Lefschetz locaux et globaux.
{\it Documents Math. Soc. Math. France} 4 (2005).

\bibitem[62]{Gr-Ra}
{\scshape L.Gruson, M.Raynaud},
Crit{\`e}res de platitude et de projectivit{\'e}.
{\it Invent. Math.} 13 (1971) pp.1--89.

\bibitem[63]{Har}
{\scshape R.Hartshorne},
Residues and duality.
{\it Springer Lect. Notes Math.} 20 (1966).

\bibitem[64]{Har2}
{\scshape R.Hartshorne},
Affine duality and cofiniteness.
{\it Invent. Math.} 9 (1969) pp.145--164.

\bibitem[65]{Har3}
{\scshape R.Hartshorne},
Generalized divisors on Gorenstein schemes.
{\it $K$-Theory} 8 (1994) pp.287--339.

\bibitem[66]{Hoch2}
{\scshape M.Hochster},
Prime ideal structure in commutative rings.
{\it Trans. Amer. Math. Soc.} 142 (1969) pp.43--60.

\bibitem[67]{Hoch}
{\scshape M.Hochster},
Rings of invariants of tori, Cohen-Macaulay
rings generated by monomials, and polytopes.
{\it Ann. of Math.} 96 (1972) pp.318--337.

\bibitem[68]{Hor}
{\scshape G.Horrocks},
Vector bundles on the punctured spectrum of a local ring.
{\it Proc. London Math. Soc.} 14  (1964) pp.689--713.

\bibitem[69]{Hu1}
{\scshape R.Huber},
Bewertungsspektrum und rigide Geometrie.
{\it Regensburger Math. Schriften} 23 (1993).

\bibitem[70]{Hu2}
{\scshape R.Huber},
Continuous valuations.
{\it Math. Zeitschr.} 212 (1993) pp.455-477.

\bibitem[71]{Hu3}
{\scshape R.Huber},
A generalization of formal schemes and rigid analytic varieties.
{\it Math. Zeithschr.} 217 (1994) pp.513--551.

\bibitem[72]{Hu-Sw}
{\scshape C. Huneke, I. Swanson},
Integral Closure of Ideals, Rings, and Modules.
{\it London Math. Soc. Lect. Notes} 336 (2006).

\bibitem[73]{Il}
{\scshape L.Illusie},
Complexe cotangent et d{\'e}formations I.
{\it Springer Lect. Notes Math.} 239 (1971).

\bibitem[74]{Il2}
{\scshape L.Illusie},
Complexe de de Rham-Witt et cohomologie crystalline.
{\it Ann. Sci. E.N.S. $4^e$ s\'erie} 12 (1979) pp.501--661.

\bibitem[75]{Il-Lz-Oz}
{\scshape L.Illusie, Y.Laszlo, F.Orgogozo},
Travaux de Gabber sur l'uniformisation locale et la cohomologie
\'etale des sch\'emas quasi-excellents.
{\it Ast\'erisque} 363--364 (2014).

\bibitem[76]{Joy}
{\scshape A.Joyal, M.Tierney},
Notes on simplicial homotopy theory.
{\it Quaderns of the CRM Barcelona} 47 (2008).

\bibitem[77]{vdKal}
{\scshape W. van der Kallen},
Descent for the K-Theory of Polynomial Rings.
{\it Math. Zeitschr.} 191 (1986) pp.405--416.

\bibitem[78]{Ka-Sch}
{\scshape M.Kashiwara, P.Schapira},
Categories and sheaves.
{\it Springer Grundl. Math. Wiss.} 332 (2006).

\bibitem[79]{Ka}
{\scshape K.Kato},
Logarithmic structures of Fontaine-Illusie.
{\it Algebraic analysis, geometry, and number theory --
Johns Hopkins Univ. Press} (1989) pp.191--224.

\bibitem[80]{Ka2}
{\scshape K.Kato},
Toric singularities.
{\it Amer. J. of Math.} 116 (1994) pp.1073-1099.

\bibitem[81]{Ked-Liu}
{\scshape K.S.Kedlaya, R.Liu},
Relative p-adic Hodge theory : foundations.
{\it Ast\'erisque} 371 (2015).

\bibitem[82]{Kel}
{\scshape J.L.Kelley},
General Topology.
{\it Springer Grad. Texts Math.} 27 (1975).

\bibitem[83]{KKMS}
{\scshape G.Kempf, F.Knudsen, D.Mumford, B.Saint-Donat},
Toroidal embeddings -- I.
{\it Springer Lect. Notes Math.} 339 (1973).

\bibitem[84]{Kn-Mu}
{\scshape F.Knudsen, D.Mumford},
The projectivity of the moduli space of stable curves --
I : Preliminaries on ``det'' and ``Div''.
{\it Math. Scand.} 39 (1976) pp.19--55.

\bibitem[85]{La}
{\scshape D.Lazard},
Autour de la platitude.
{\it Bull. Soc. Math. France} 97 (1969) pp.81--128.

\bibitem[86]{Lip}
{\scshape J.Lipman},
Rational singularities with applications to algebraic
surfaces and unique factorization.
{\it Publ.Math. IHES} 36 (1969) pp.195--279.

\bibitem[87]{Mas}
{\scshape W.S.Massey},
A basic course in algebraic topology.
{\it Springer GTM} 127 (1991).

\bibitem[88]{Mat-old}
{\scshape H.Matsumura},
Commutative algebra -- Second edition.
{\em Benjamin/Cummings Math. Lect. Notes} (1980).

\bibitem[89]{Mat}
{\scshape H.Matsumura},
Commutative ring theory.
{\it Cambridge Univ. Press} (1986).

\bibitem[90]{Mo}
{\scshape S.Mochizuki},
Extending families of curves over log regular schemes.
{\it J. Reine Angew. Math.} 511 (1999) pp.43--71.

\bibitem[91]{Mu}
{\scshape D.Mumford},
Lectures on curves on algebraic surfaces.
{\it Princeton Univ. Press} (1966).

\bibitem[92]{Na}
{\scshape M.Nagata}, Local rings.
{\it Interscience Publ.} (1962).

\bibitem[93]{Na2}
{\scshape M.Nagata}, Finitely generated rings over a valuation
ring {\it J.Math. Kyoto Univ.} (1966) pp.163--169.

\bibitem[94]{Niz}
{\scshape W.Niziol},
Toric singularities : log-blow-ups and global resolutions.
{\it J. of Alg. Geom.} 15 (2006) pp.1--29.

\bibitem[95]{Oh-Pe}
{\scshape J.Ohm, R.L.Pendelton},
Rings with noetherian spectrum.
{\it Duke Math. J.} 35 (1968) pp.631--639.

\bibitem[96]{Ol}
{\scshape J.-P.Olivier},
Going up along absolutely flat morphisms.
{\it J.Pure Appl. Algebra} 30 (1983) pp.47--59.

\bibitem[97]{Qu0}
{\scshape D.Quillen},
Homology of commutative rings, {\it unpublished notes} MIT (1968).  

\bibitem[98]{Qu}
{\scshape D.Quillen},
Higher algebraic $K$-theory -- I.
{\it Springer Lect. Notes Math.} 341 (1973) pp.85--147.

\bibitem[99]{Ray}
{\scshape M.Raynaud},
Anneaux locaux hens{\'e}liens.
{\it Springer Lect. Notes Math.} 169 (1970).

\bibitem[100]{Ro}
{\scshape N.Roby},
Lois polyn\^omes et lois formelles en th\'eorie des modules.
{\it Ann. Sci. E.N.S. $3^e$ s\'erie}  80 (1963) pp.213--348.

\bibitem[101]{Z-S}
{\scshape P.Samuel, O.Zariski},
Commutative algebra, Vol.II.
{\it Springer Grad. Texts Math.} 29 (1975).

\bibitem[102]{Sch}
{\scshape C.Scheiderer},
Quasi-augmented simplicial spaces, with an application
to cohomological dimension.
{\it J.Pure Appl. Algebra} 81 (1992) pp.293--311.

\bibitem[103]{Scho}
{\scshape P.Scholze},
Perfectoid spaces.
{\it Publ. Math. IHES} 116 (2012) pp.245--313.

\bibitem[104]{Scho2}
{\scshape P.Scholze},
$p$-adic Hodge theory for rigid-analytic spaces.
{\it Forum of Mathematics, Pi} 1 (2013) e1; corrigendum 4 (2016) e6.

\bibitem[105]{Sp}
{\scshape N.Spaltenstein},
Resolutions of unbounded complexes.
{\it Comp. Math.} 65 (1988) pp.121--154.

\bibitem[106]{Sta}
{\scshape R.P.Stanley},
Hilbert functions of graded algebras.
{\it Adv. in Math.} 28 (1978) pp.57--83.

\bibitem[107]{Sw}
{\scshape R.Swan},
On seminormality.
{\it J. Alg.} 67 (1980) pp.210--229.

\bibitem[108]{Tem}
{\scshape M.Temkin},
Relative Riemann-Zariski spaces.
{\it Israel J. of Math.} 185 (2011) pp.1-42.

\bibitem[109]{Tsu}
{\scshape T.Tsuji},
Saturated morphisms of logarithmic schemes.
{\it Tunisian J. of Math.} 1 (2019) pp.185--220.

\bibitem[110]{Yeku}
{\scshape R.Vyas, A.Yekutieli},
Weak proregularity, weak stability, and the non-commutative
MGM equivalence.
{\it Preprint on} {\ttfamily arXiv:1608.03543} (2017).

\bibitem[111]{Yan}
{\scshape H.Yanagihara},
Some results on weakly normal ring extensions.
{\it J.Math. Soc. Japan} 35 pp.649--661 (1983).
  
\bibitem[112]{We}
{\scshape C.Weibel},
An introduction to homological algebra.
{\it Cambridge Univ. Press} (1994).

\end{thebibliography}
\end{document}